\documentclass[a4paper, 10pt, reqno]{amsbook}
\usepackage{BOONDOX-cal}
\usepackage{fontawesome}
\usepackage{wrapfig}
\usepackage{epigraph} 
\usepackage{amssymb}
\usepackage{environ}
\usepackage{amsmath}
\usepackage[mathscr]{eucal}
\usepackage{mathrsfs}
\usepackage{xspace} 
\usepackage{adjustbox}
\usepackage{tabu}
\usepackage{thmtools}
\usepackage[utf8]{inputenc}
\usepackage[english]{babel}
\usepackage[all,2cell]{xy}
\usepackage{xcolor}
\usepackage{graphicx}
\usepackage{wasysym} 
\usepackage{hyperref}
\usepackage{pdflscape}
\usepackage{caption}
\usepackage{multicol}
\usepackage{longtable}

\numberwithin{section}{chapter}
\numberwithin{subsection}{section}
\numberwithin{equation}{section}
\numberwithin{figure}{chapter}

\makeatletter
\newcommand{\leqnomode}{\tagsleft@true}
\newcommand{\reqnomode}{\tagsleft@false}

\newsavebox{\measure@tikzpicture}
\NewEnviron{scaletikzpicturetowidth}[1]{%
  \def\tikz@width{#1}%
  \def\tikzscale{1}\begin{lrbox}{\measure@tikzpicture}%
  \BODY
  \end{lrbox}%
  \pgfmathparse{#1/\wd\measure@tikzpicture}%
  \edef\tikzscale{\pgfmathresult}%
  \BODY
}

\makeatother

\usepackage{array}
\newcolumntype{L}[1]{>{\raggedright\let\newline\\\arraybackslash\hspace{0pt}}m{#1}}
\newcolumntype{C}[1]{>{\centering\let\newline\\\arraybackslash\hspace{0pt}}m{#1}}
\newcolumntype{R}[1]{>{\raggedleft\let\newline\\\arraybackslash\hspace{0pt}}m{#1}}

\definecolor{sg}{HTML}{df78ef}
\definecolor{sg1}{HTML}{ab47bc}
\definecolor{sg2}{HTML}{790e8b}

\usepackage{tikz}
\usetikzlibrary{mindmap,backgrounds, calc}
\usetikzlibrary{matrix,chains,scopes,positioning,arrows,fit}
\usetikzlibrary{positioning, shapes, patterns}
\usetikzlibrary{automata} 
\usetikzlibrary{shapes.geometric, arrows,arrows.meta}
\usetikzlibrary{positioning,decorations.markings}
\usepackage{hyperref}
\usepackage[acronym]{glossaries}

\usepackage{imakeidx}
\makeindex

\usepackage{stmaryrd} 
\UseAllTwocells

\makeatletter
\newcommand*\bigcdot{\mathpalette\bigcdot@{.5}}
\newcommand*\bigcdot@[2]{\mathbin{\vcenter{\hbox{\scalebox{#2}{$\m@th#1\bullet$}}}}}
\makeatother

\newcommand{\at}[1]%
            {\ensuremath{\protect\underline{\mathbf{#1}}}} 

\newcommand{\op}[1]{\ensuremath{\operatorname{#1}}}        

\newcommand{\h}[1][]                                       
 {\ifthenelse{\boolean{mmode}}%
  {$\mathrm{h}$}%
  {h\nobreakdash#1\hspace{0pt}}}

\newcommand{\id}{\mathrm{id}}        
\DeclareMathOperator{\Hom}{Hom}    
\newcommand{\comp}{\circ}          
\newcommand{\adcomp}%
  {\overset{\operatorname{ad}}{\comp}} 
\newcommand{\funcomp}%
  {\overset{\operatorname{fn}}{\comp}}

\newcommand{\sccat}
{\mathbin{\kern-1pt\raisebox{6pt}{.}\kern-5pt
\downarrow\kern-5pt\raisebox{6pt}{.}\kern-1pt}}

\newcommand{\parrow}[1]
   {\underset{{\displaystyle \raisebox{5pt}%
   {$\longleftarrow$}}}{\op{#1}}{\,}}
\newcommand{\iarrow}[1]
   {\underset{{\displaystyle \raisebox{5pt}%
   {$\longrightarrow$}}}{\op{#1}}{\,}}

\newcommand{\uadj}{\top}           

\DeclareMathOperator{\Sub}{Sub}    
\DeclareMathOperator{\Rel}{Rel}    
\DeclareMathOperator{\Fnc}{Fnc}    
\newcommand{\rest}%
{\mathnormal{\restriction}}        

\newcommand{\nin}{\not\in}         
\newcommand{\bprod}{\times}        
\newcommand{\iso}{\cong}           
\newcommand{\function}[4]{
            \begin{array}{@{\:}c@{\:}c@{\:}l}
                   #1 &\mor& #2 \\
                   #3 &\longmapsto& #4
            \end{array} }
\newcommand{\nfunction}[4]
    {\left\{
     \function{#1}{#2}{#3}{#4}
     \right. }

\newcommand{\bb}[1]{\ensuremath
 {\lvert #1 \rvert}}
\DeclareMathOperator{\Sg}{Sg}      
\DeclareMathOperator{\Ker}{Ker}    
\newcommand{\pr}{\mathrm{pr}}        
\newcommand{\vs}[1]{\mathbin{\downarrow}#1}
\DeclareMathOperator{\G}{G}        
\DeclareMathOperator{\bconcat}
            {\curlywedge}

\newcommand{\concat}
  {\ensuremath{\text
  {\Large $\curlywedge$}}}
\newcommand{\ext}[1]
  {\ensuremath{#1^{\sharp}}}

\DeclareMathOperator{\T}{T}        

\newcommand{\ol}{\overline}

\newcommand{\brel}{\ensuremath{\xymatrix{{}\arity@{{*}{-}{*}}[r] & {}}}}

\newcommand{\nseq}[3]{\xymatrix@1@C=16pt{#1 \arity@{>}[r]_-{\scriptscriptstyle{#2}} & #3 }}

\NoCompileMatrices              
\xymatrixrowsep = {8ex}         
\xymatrixcolsep = {10ex}        

\newdir{<= }{:a(+25)\dir{-}*:a(-25)\dir{-}*!/:a(85)+1.8pt/:a(-25)\dir{-}}
\newdir{ = }{*@{=}}
\newdir{+>}{@{|}*@{-}!/-6.5pt/@{>}*!/-13pt/{}}
\newdir{ +}{{}*!/-5pt/@{+}}
\newdir{ >}{{}*!/-10pt/@{>}}
\newdir{.-}{{}*:a(-90)h!/2.1pt/{.}}
\newdir{.->}{:a(-90)h!/2.5pt/{.}*!/4pt/\dir{-}*!/0pt/\dir{-}*%
!/-2pt/\dir{-}*!/-8pt/\dir{>}}
\newdir{-.->}%
{:a(-90)h!/2pt/{.}*!/8pt/\dir{-}*!/4pt/\dir{-}*!/0pt/\dir{-}*%
!/-5pt/\dir{-}*!/-11pt/\dir{>}}
\newdir{~>}{!/4.6pt/\dir{~}*!/1pt/\dir{-}*!/-5pt/\dir{>}}
\newdir{~~>}{!/7pt/\dir{~}*!/0pt/\dir{~}*!/-5pt/\dir{>}}
\newdir{=~>}{*!/1pt/@2{~}*!/6pt/{}*!/-6pt/@2{>}}
\newdir{<~~>}{!/21pt/\dir{<}*!/21pt/\dir{-}!/12pt/\dir{~}!/5pt/\dir{~}*!/1pt/\dir{-}*!/-3pt/\dir{>}}
\newdir{=>}{!/5pt/\dir{=}!/2.5pt/\dir{=}*!/-5pt/\dir2{>}*!/7pt/\dir{ }}
\newdir{==>}{!/-2pt/\dir{=}!/-6pt/\dir{=}%
 !/-10pt/\dir{=}*!/-18pt/\dir2{>}}
\newdir{<==>}{!\dir2{<}!/-2pt/\dir{=}!/-6pt/\dir{=}%
 !/-10pt/\dir{=}*!/-17pt/\dir2{>}}
\newdir{< }{{}*!/12pt/\dir{<}}
\newdir{-o}{!/1.3pt/{}*!/0pt/{}@{o}*!/-5pt/{}}

\newsavebox{\xymor}  
\newsavebox{\xymon}  
\newsavebox{\xyepi}  
\newsavebox{\xytn}   
\newsavebox{\xyrel}  
\newsavebox{\xycel}  
\newsavebox{\xymdf}  
\newsavebox{\xyumor} 
\newsavebox{\xydmor} 
\newsavebox{\xyomor} 
\newsavebox{\xyemor} 

\newcommand{\xynode}{\makebox[0ex]{}}
\savebox{\xymor}{\ensuremath{%
\xymatrix@1@C=19pt{\xynode \ar@{>}[r] & \xynode }}}
\savebox{\xymon}{\ensuremath{%
\xymatrix@1@C=19pt{\xynode \ar@{{ +}{-}{>}}[r] & \xynode }}}
\savebox{\xyepi}{\ensuremath{%
\xymatrix@1@C=19pt{\xynode \ar@{{}{-}{+>}}[r] & \xynode }}}
\savebox{\xytn}{\ensuremath{%
\xymatrix@1@C=19pt{\xynode \ar[r]|(.44){\object@{.-}} & \xynode
}}}
\savebox{\xyrel}{\ensuremath{%
\xymatrix@1@C=19pt{\xynode \ar@{{}{-}{-o}}[r] & \xynode }}}
\savebox{\xycel}{\ensuremath{%
\xymatrix@1@C=19pt{\xynode \ar@{=>}[r] & \xynode }}}
\savebox{\xymdf}{\ensuremath{%
\xymatrix@1@C=16pt{\xynode \ar@{}[r]|{\dir{~>}} & \xynode}}}
\savebox{\xyumor}{\ensuremath{%
\xymatrix@1@C=19pt{\xynode \ar@{{}{-}^{>}}[r] & \xynode }}}
\savebox{\xydmor}{\ensuremath{%
\xymatrix@1@C=19pt{\xynode \ar@{{}{-}_{>}}[r] & \xynode }}}
\savebox{\xyomor}{\ensuremath{%
\xymatrix@1@C=19pt{\xynode \ar@{{}{-}^{< }}[r] & \xynode }}}
\savebox{\xyemor}{\ensuremath{%
\xymatrix@1@C=19pt{\xynode \ar@{{ >}{-}{>}}[r] & \xynode }}}

\newcommand{\mor}{\usebox{\xymor}}    
\newcommand{\mon}{\usebox{\xymon}}    
\newcommand{\epi}{\usebox{\xyepi}}    
\newcommand{\cel}{\usebox{\xycel}}    
\newcommand{\dmor}{\usebox{\xydmor}}  

\newcommand{\functor}[9]{
 \xymatrix{
    #4 \save[]+<0ex,5ex>*+{#1}="1"  \restore
      \arity[d]_{#6}  \arity@{}[rd]|{\longmapsto}
  & #5 \save[]+<0ex,5ex>*+{#3}="3"  \restore
      \arity[d]^{#7}
  \\
   #8 & #9 \arity "1";"3"^-{#2} } }
\newcommand{\functornd}[9]{
 \xymatrix{
    #4 \save[]+<0ex,5ex>*+{#1}="1"  \restore
      \arity[d]_{#6}  \arity@{}[rd]|{\longmapsto}
  & #5 \save[]+<0ex,5ex>*+{#3}="3"  \restore
  \\
   #8 & #9 \arity[u]_{#7} \arity "1";"3"^-{#2} } }
\newcommand{\functordn}[9]{
 \xymatrix{
    #4 \save[]+<0ex,5ex>*+{#1}="1"  \restore
       \arity@{}[rd]|{\longmapsto}
  & #5 \save[]+<0ex,5ex>*+{#3}="3"  \restore
      \arity[d]^{#7}
  \\
   #8  \arity[u]^{#6}  & #9 \arity "1";"3"^-{#2} } }

\newcommand{\larr}{->}
\newcommand{\rarr}{->}

\newcommand{\xfunctor}[9]{
 \xymatrix{
    #4 \save[]+<0ex,5ex>*+{#1}="1"  \restore
      \ifthenelse{\equal{\larr}{->}}{\arity[d]_{#6}}{}
      \ifthenelse{\equal{\larr}{<-}}{\arity[d];[]^{#6}}{}
      \ifthenelse{\equal{\larr}{-<}}{\arity@{< }[d]_{#6}}{}
      \arity@{}[rd]|{\longmapsto}
  & #5 \save[]+<0ex,5ex>*+{#3}="3"  \restore
      \ifthenelse{\equal{\rarr}{->}}{\arity[d]^{#7}}{}
      \ifthenelse{\equal{\rarr}{<-}}{\arity[d];[]_{#7}}{}
      \ifthenelse{\equal{\rarr}{-<}}{\arity@{< }[d]^{#7}}{}
  \\
   #8 & #9 \arity "1";"3"^-{#2} } }

\theoremstyle{plain}
\newtheorem{theorem}{Theorem}[section]
\newtheorem{proposition}[theorem]{Proposition}

\newtheorem{corollary}[theorem]{Corollary}
\newtheorem{lemma}[theorem]{Lemma}
\newtheorem{remark}[theorem]{Remark}
\theoremstyle{definition}
\newtheorem{definition}[theorem]{Definition}
\newtheorem{example}[theorem]{Example}

\newtheorem{assumption}{\bf Assumption}

\newtheorem*{convention}{Convention}

\newcommand{\arity}{\mathsf{ar}}

\newtheorem{claim}[theorem]{Claim}

\usepackage[LGR,T1]{fontenc} 
\usepackage{etoolbox} 
\DeclareSymbolFont{sfletters}{OML}{cmbrm}{m}{it} 
\DeclareMathSymbol{\w}{\mathord}{sfletters}{'041}      
\begin{document}
\title[Higher-order rewriting systems]{From higher-order rewriting systems to higher-order categorial algebras and higher-order Curry-Howard isomorphisms}

\author[Climent]{J. Climent Vidal}
\address{Universitat de Val\`{e}ncia\\
         Departament de L\`{o}gica i Filosofia de la Ci\`{e}ncia\\
         Av. Blasco Ib\'{a}\~{n}ez, 30-$7^{\mathrm{a}}$, 46010 Val\`{e}ncia, Spain}
\email{Juan.B.Climent@uv.es}
\author[Cosme]{E. Cosme Ll\'{o}pez}
\address{Universitat de Val\`{e}ncia\\
         Departament de Matem\`{a}tiques\\
         Dr. Moliner, 50, 46100 Burjassot, Val\`{e}ncia, Spain}
\email{Enric.Cosme@uv.es}

\date{January 2026}

\begin{abstract}
This ongoing project aims to comprehensively explore the study of higher-order rewriting systems and higher-order categorial algebras. The main objective is to develop a solid theoretical framework that serves as a foundation for future research and practical applications. To achieve this, higher-order Curry-Howard type isomorphisms are introduced. The current version fully describes the construction of free categorial algebras associated with rewriting systems of orders 0, 1, and 2.  The following text is open to discussions within the scientific community.  All improvements and suggestions are welcome.

\end{abstract}


\begin{titlepage}

\centering

\vspace{1cm}
{\huge\bfseries From higher-order rewriting systems to higher-order categorial algebras and higher-order Curry-Howard isomorphisms\par}

\vspace{3cm}

{\Large \textbf{Juan Climent Vidal}\par}
\vspace{.25cm}
{\itshape 	Universitat de Val\`{e}ncia\\
			Departament de L\`{o}gica i Filosofia de la Ci\`{e}ncia\par}
{\small Juan.B.Climent@uv.es\par}
         	
\vspace{.75cm}

{\Large \textbf{Enric Cosme Ll\'{o}pez}\par}
\vspace{.25cm}
{\itshape 	Universitat de Val\`{e}ncia\\
			Departament de Matem\`{a}tiques\par}
{\small Enric.Cosme@uv.es (\faEnvelopeO)\par}
{\,\\}
{\itshape 	Nantong University\\
			School of Mathematics and Statistics\par}

\vspace{.75cm}

{\Large \textbf{Ra\'{u}l Ruiz Mora}\par}
\vspace{.25cm}
{\itshape 	Universitat de Val\`{e}ncia\\
			Departament de Matem\`{a}tiques\par}
{\small Raul.Ruiz-Mora@uv.es\par}

\vfill

{\small January 2026\par}

\vspace{1cm}


\end{titlepage}

\newpage
\thispagestyle{empty} 
{\qquad}

\newpage
\thispagestyle{empty} 
{\qquad}

\vspace{1.75cm}
\begin{epigraph}
{
\noindent
Or la possibilit\'{e} de la traduction implique l'existence d'un invariant. Traduire, c'est pr\'{e}cisément d\'{e}gager cet invariant.
\\
\emph{[Now the possibility of translation implies the existence of an invariant. To translate is precisely to disengage this invariant.]}
}{
Henri Poincar\'{e},\\
\emph{La Valeur de la Science,} 1905.
}
\end{epigraph}

\vspace{.75cm}
\begin{epigraph}
{
\noindent
Caminante, no hay camino, \\
se hace camino al andar.\\
Al andar se hace el camino,\\
y al volver la vista atr\'{a}s\\
se ve la senda que nunca\\
se ha de volver a pisar.
\\
\emph{[Traveler, there is no path,\\ the path is made by walking.\\By walking the path is made\\and when you look back\\you'll see a road\\never to be trodden again.]}
}{
Antonio Machado,\\
\emph{Proverbios y cantares,} 1912.
}
\end{epigraph}

\vspace{.75cm}
\begin{epigraph}
{
\noindent
El clar cam\'{i}, el pregon idioma, \\
    un alfabet fosforescent de pedres, \\
    un alfabet sempre amb la clau al pany.
\\
\emph{[The clear path, the deep language,\\a phosphorescent alphabet of stones,\\an alphabet always with the key in the lock.]}
}{
Vicent Andr\'{e}s Estell\'{e}s,\\
\emph{X\`{a}tiva,} 1980.
}
\end{epigraph}

\newpage
\thispagestyle{empty} 
{\qquad}

\newpage
\thispagestyle{empty} 

\begin{center}
{\large\bfseries Summary\par}
\end{center}
       
\vspace{.5cm} 
This ongoing project aims to define and investigate, from the standpoint of category theory, order theory, and universal algebra, the notions of higher-order many-sorted rewriting system and of higher-order many-sorted categorial algebra and their relationships, via the higher-order Curry-Howard isomorphisms. The work carried out effectively so far is divided into three parts. This summary will outline the content of each part.

In \textsf{Part 0}, for the convenience of the reader, we gather together those facts about many-sorted sets, many-sorted (total and partial) algebras, category theory, order theory, and higher-order categories, which will be used extensively in the rest of the paper. 
Readers who are familiar with the just mentioned fields may wish to skip this part and return to it as required.

In \textsf{Part 1}, we define the many-sorted set of paths associated with a many-sorted rewriting system and equip it with a structure of many-sorted partial algebra  and a structure of category, and to its coproduct with a structure of Artinian ordered set. Next, we consider an extension of the signature associated with the many-sorted rewriting system, and we associate each path with a term in the extended signature. This constitutes a Curry-Howard type mapping. After that we prove that the quotient of the many-sorted set of paths by the kernel of the Curry-Howard mapping is equipped with a structure of many-sorted partial algebra and a structure of category, and that its coproduct is equipped with a structure of Artinian ordered set. Following this we identify a subquotient of the free term algebra in the extended signature that is isomorphic to the algebraic, categorical, and ordered structures on the quotient of paths. Additionally, we prove that these two structures are isomorphic to the free partial algebra on paths in a variety of partial algebras for the extended signature.  This will be a Curry-Howard type isomorphism.

In \textsf{Part 2}, following the homotopic dictum, we define the many-sorted set of second-order paths, i.e., paths on paths, associated with a second-order rewriting system and delve into the structures of partial algebra, category, and Artinian order that can be defined upon it. Next, we consider a second-order extension of the signature associated with the second-order rewriting system, and we associate each second-order path with a term in the extended second-order signature. This constitutes a second-order Curry-Howard type mapping. It is proven that a refined quotient of the many-sorted set of second-order paths by the kernel of the second-order Curry-Howard mapping admits a structure of a partial algebra, $2$-category, and Artinian preordered set. Next, we identify a subquotient of the free algebra of terms in the extended second-order signature that is isomorphic to the algebraic, $2$-categorical, and ordered structures on the quotient of second-order paths. Additionally, we prove that these two structures are isomorphic to the free partial algebra generated by second-order paths in a variety of partial algebras for the second-order extended signature.   This would be a second-order Curry-Howard type isomorphism.

In \textsf{Part 3}, we develop a tiered definition for morphisms between rewriting systems. At layer $0$, which encompasses the underlying many-sorted algebraic basis—sorts, operations, and variables—we employ the concept of a derivor to formalize the notion of a zeroth-order morphism. We introduce the extension of this definition to terms, subsequently constructing the category of zeroth-order rewriting systems. This framework is extended to layer $1$ by incorporating rewriting rules; here, a first-order morphism is defined via the preceding zeroth-order morphism coupled with a mapping that assigns the rewriting rules in the domain to paths within the codomain. We show that this mapping can be extended to paths and that  this extension is compatible with path classes, leading to the construction of the category of first-order rewriting systems. Finally, this inductive process is iterated at layer $2$ to account for second-order rewriting rules. A second-order morphism is characterized by a first-order morphism and a mapping that assigns  second-order rewriting rules in the domain to second-order  paths in the  codomain. By extending this to second-order paths and ensuring compatibility with second-order path classes, we construct the category of second-order rewriting systems.

The ultimate goal, to be developed in future versions of this work, is to define and investigate the category of towers, whose objects will consist of families, indexed by $\mathbb{N}$, of higher-order many-sorted rewriting systems and of higher-order many-sorted categorial algebras, including higher-order Curry-Howard type results for the latter, together with an additional structure that intertwines such $\mathbb{N}$-families; and whose morphism from a tower to another will be families, indexed by $\mathbb{N}$,  of morphisms between its higher-order many-sorted rewriting systems and of higher-order many-sorted categorial algebras compatible with their structures. 
All feedback is appreciated.

\thispagestyle{empty} 
{\qquad}
\vfill

{\footnotesize
  \par}

\clearpage
\setcounter{page}{-1}   

\frontmatter
\setcounter{tocdepth}{1}
\tableofcontents 

\mainmatter
\part{Preliminaries}
\chapter{Introduction}\label{S0A}


The theory of single-sorted term rewriting systems (which for words has its origins in the works of Thue~\cite{Th14}, Dehn~\cite{MD11} and Post~\cite{Post43,Post65}), changing one term to another according to certain rewrite rules or productions, 
is a fundamental field
within computer science. Moreover, it has proved to be effective in solving decision problems in a wide range of areas, including mathematical logic, mathematics and physics (such as, e.g., the word problem for semigroups, unsolvability theorems for logics, Hilbert's tenth problem and undecidability in classical mechanics). Briefly stated, it could be said that rewriting is the root of all computational processes. 

On the other hand, the classical Curry-Howard correspondence explains the direct relationship between computer programs and mathematical proofs. More precisely, Curry, in~\cite{CF58}, was the first to acknowledge the formal analogy between his combinatory logic and the axioms of a Hilbert--type deduction system for the positive implicational propositional logic. Later on, Howard, in 1969, but published in~\cite{H80}, observed the same formal analogy between Church's $\lambda$-calculus and the proof rules of a Gentzen's system of natural deduction for intuitionistic propositional logic. The Curry-Howard correspondence, also called correspondence proofs-terms or correspondence proofs-programs, assigns to each proof in the intuitionistic logic a term in Curry's combinatory logic or in Church's $\lambda$-calculus. In other words, the Curry-Howard correspondence consists of the clear-sighted observation that two seemingly unrelated families of formalisms---namely, systems of formal deduction, on the one hand, and models of computation, on the other---are, essentially, the same kind of mathematical object.

What we present in this work is the preliminary development of a theory aimed at defining the notions of higher-order many-sorted rewriting systems and higher-order many-sorted categorial algebras and investigating the relationship between them  through higher-order many-sorted Curry-Howard isomorphisms (please see below for a description of the work we intend to develop in the near future). Our frequent use of the qualifier Curry-Howard in this paper is due to the fact that we have been able to represent paths, that is, syntactic derivations between terms for a rewriting system (which have a proof-theoretical flavor) as terms of an algebra relative to a signature associated with the rewriting system. Furthermore, following the homotopic dictum, we have also found a way to iterate this process so that higher-order paths correspond to terms in higher-order signatures. On the other hand, paraphrasing a sentence of N. Gambino and A. Joyal in~\cite{GJ17}, we can say that dropping the restriction of working with single-sorted algebras and single-sorted term rewriting systems and allowing higher-order many-sorted (partial and total) algebras and higher-order many-sorted term rewriting systems not only does not create problems---except those notational ones that inevitably arise when raising the order---, but in fact brings to light new mathematical structures that were hidden in the world of single-sorted term rewriting systems. 

Let us explain in what sense we use the term ``higher'' in this work.
For a set of sorts $S$ and an $S$-sorted signature $\Sigma$, our definition of the category $\mathsf{nCat}^{S}\mathsf{Alg}(\Sigma)$, of $S$-sorted strict $n$-categorial
$\Sigma$-algebras will be such that the forgetful functor from it to $\mathsf{Set}^{S}$, the category of $S$-sorted sets, factorizes through the category $\mathsf{nCat}^{S}$, of $S$-sorted strict $n$-categories. It is precisely in this sense in which the term ``higher'' is used in ``higher-order categorial algebras'', and, derivatively, in ``higher-order rewriting systems'', because, as will be seen below, to every higher-order rewriting system one can naturally associate a (labelled) higher-order categorial algebra. Let us note that, for every $n\in \mathbb{N}$, the category $\mathsf{nCat}$, of strict $n$-categories, is a particular case of the categories of type $\mathsf{nCat}^{S}\mathsf{Alg}(\Sigma)$ (take as $S$ a final set, e.g., the ordinal $1 = \{0\}$, and as $\Sigma$ the initial signature).  

In this work, from very basic theories and facts (such as category theory, the theory of ordered sets, many-sorted total or partial algebra and the construction of the free completion of a many-sorted partial algebra), we obtain, for $n=0,1,2$, pairs of many-sorted $n$-categorial algebras that are isomorphic to each other by means of many-sorted $n$-th order Curry-Howard type isomorphisms.  In our case, the different Curry-Howard isomorphisms  identify equivalence classes of $n$-th order paths and equivalence classes of terms, what we have called $n$-th order path term classes. In this regard, it seems pertinent to recall what P. Roquette in~\cite{pr18}, on p.~1, wrote:  ``Sometimes a result appears to be better understood if it is generalized and freed from unnecessary assumptions, or if it is embedded into a general theory which opens analogies to other fields of mathematics. Also, in order to make further progress possible it is often convenient and sometimes necessary to develop a framework, conceptual and notational, in which the known results become trivial and almost self-evident, at least from a formalistic point of view.''  Further clarification on these concepts and others will be provided in the subsequent chapters and thoroughly discussed following the introduction.

\section{Our work}
The paper is divided into four parts, and below we provide a more detailed explanation of the notions and constructs introduced in it. In addition, we provide the guidelines for what will be the natural continuation of the work done so far.



\subsection{Preliminaries}

\textsf{Part~0} starts by reviewing fundamental notions and facts on many-sorted sets and algebras. These include the constructive description of the subalgebra and congruence generating operators, as well as that of the free many-sorted algebra on a many-sorted set. Furthermore, the many-sorted translation operators and their generalizations are defined. Afterwards, a detailed summary of the most important results on partial many-sorted algebras is provided. This includes the discussion of several types of homomorphism, subalgebras, and congruences. We also introduce the free completions associated to partial many-sorted algebras as well as the varieties of partial many-sorted algebras and the description of the associated free solutions.
Then the notions of adjunction, universal morphism from an object to a functor, map from an adjunction to another, the categories of split indexed category, split fibred category and the construction of Grothendieck. 
Moreover, after recalling the notions of (single-sorted) strict $n$-category and (single-sorted) strict $\omega$-category and their many-sorted versions, a proof of their respective equivalences is given. 
Finally,  the notions of many-sorted strict $n$-categorical $\Sigma$-algebra and strict $n$-categorical $\Sigma$-homomorphism between them are defined.
Readers who are familiar with the above may wish to skip this part and return to it as required.

\subsection{First-order structures associated to a rewriting system}
\textsf{Part~1} begins by defining a many-sorted term rewriting system $ \boldsymbol{\mathcal{A}}^{(1)}$, written simply as $\boldsymbol{\mathcal{A}}$, as an ordered pair $\boldsymbol{\mathcal{A}} = (\boldsymbol{\mathcal{A}}^{(0)},\mathcal{A})$, where $\boldsymbol{\mathcal{A}}^{(0)} = (S,\Sigma,X)$ is a zeroth order $S$-sorted rewriting system, i.e., $S$ is a set of sorts, $\Sigma$ an $S$-sorted signature, $X$ an $S$-sorted set, and $\mathcal{A}$ a subset of $(\mathrm{T}_{\Sigma}(X)^{2}_{s})_{s\in S}$, the $S$-sorted set of the rewrite rules  with variables in $X$, where $\mathrm{T}_{\Sigma}(X)$ is the underlying $S$-sorted set of $\mathbf{T}_{\Sigma}(X)$, the free $\Sigma$-algebra on $X$. Then we define the notion of path in $\boldsymbol{\mathcal{A}}$ from a term to another. Following this we show that paths have a recursive decomposition and we call this algorithm the path extraction algorithm. Next, we define the $S$-sorted set $\mathrm{Pth}_{\boldsymbol{\mathcal{A}}}$, of paths in 
$\boldsymbol{\mathcal{A}}$, and equip it with a structure of $\Sigma$-algebra. Moreover, we equip $\coprod \mathrm{Pth}_{\boldsymbol{\mathcal{A}}}$, the coproduct of $\mathrm{Pth}_{\boldsymbol{\mathcal{A}}}$, with a structure of Artinian order. After this, from $\boldsymbol{\mathcal{A}}$, we obtain an $S$-sorted signature $\Sigma^{\boldsymbol{\mathcal{A}}}$, which is an extension of the $S$-sorted signature $\Sigma$, and equip $\mathrm{Pth}_{\boldsymbol{\mathcal{A}}}$ with a structure of $S$-sorted partial $\Sigma^{\boldsymbol{\mathcal{A}}}$-algebra, that we denote by $\mathbf{Pth}_{\boldsymbol{\mathcal{A}}}$. Next, we define, by Artinian recursion, an $S$-sorted mapping $\mathrm{CH}^{(1)}$, the $S$-sorted Curry-Howard mapping, from the $S$-sorted set 
$\mathrm{Pth}_{\boldsymbol{\mathcal{A}}}$ to $\mathrm{T}_{\Sigma^{\boldsymbol{\mathcal{A}}}}(X)$, the underlying $S$-sorted set of the free $S$-sorted 
$\Sigma^{\boldsymbol{\mathcal{A}}}$-algebra $\mathbf{T}_{\Sigma^{\boldsymbol{\mathcal{A}}}}(X)$ on $X$. We  prove that the Curry-Howard mapping is not a $\Sigma^{\boldsymbol{\mathcal{A}}}$-homomorphism from $\mathbf{Pth}_{\boldsymbol{\mathcal{A}}}$ to $\mathbf{T}_{\Sigma^{\boldsymbol{\mathcal{A}}}}(X)$. 
However, by a suitable passage to the quotient, both in the domain and in the codomain of $\mathrm{CH}^{(1)}$, we obtain three pairs of first order structures, of $S$-sorted partial $\Sigma^{\boldsymbol{\mathcal{A}}}$-algebras, of $S$-sorted categorial $\Sigma$-algebras and of Artinian orders, respectively, such that a derived $S$-sorted mapping $\mathrm{CH}^{[1]}$ of $\mathrm{CH}^{(1)}$ induces isomorphisms between the first-order structures of each one of the above pairs.

To prove that $\mathrm{CH}^{[1]}$ induces the aforementioned isomorphisms, we proceed as follows. In the first place, after showing that $\mathrm{Ker}(\mathrm{CH}^{(1)})$, the kernel of $\mathrm{CH}^{(1)}$, is a closed $\Sigma^{\boldsymbol{\mathcal{A}}}$-congruence on $\mathbf{Pth}_{\boldsymbol{\mathcal{A}}}$, we define the quotient $[\mathrm{Pth}_{\boldsymbol{\mathcal{A}}}]$ of $\mathrm{Pth}_{\boldsymbol{\mathcal{A}}}$ by $\mathrm{Ker}(\mathrm{CH}^{(1)})$ and prove that is equipped with a structure of $S$-sorted partial $\Sigma^{\boldsymbol{\mathcal{A}}}$-algebra, a structure of $S$-sorted categorial $\Sigma$-algebra and, its coproduct, with a structure of Artinian order. We denote by $[\mathbf{Pth}_{\boldsymbol{\mathcal{A}}}]$, $[\mathsf{Pth}_{\boldsymbol{\mathcal{A}}}]$ and $(\coprod [\mathrm{Pth}_{\boldsymbol{\mathcal{A}}}],\leq_{[\mathbf{Pth}_{\boldsymbol{\mathcal{A}}}]})$, respectively, the corresponding $S$-sorted partial $\Sigma^{\boldsymbol{\mathcal{A}}}$-algebra, $S$-sorted categorial $\Sigma$-algebra and Artinian ordered set. 

Afterwards, we define a congruence $\Theta^{[1]}$ on $\mathbf{T}_{\Sigma^{\boldsymbol{\mathcal{A}}}}(X)$ and, from it, we obtain the $S$-sorted subset $\mathrm{PT}_{\boldsymbol{\mathcal{A}}}$, of path terms, of $\mathrm{T}_{\Sigma^{\boldsymbol{\mathcal{A}}}}(X)$ as the $\Theta^{[1]}$-saturation of $\mathrm{CH}^{(1)}[\mathrm{Pth}_{\boldsymbol{\mathcal{A}}}]$. Then we define $[\mathrm{PT}_{\boldsymbol{\mathcal{A}}}]$, the $S$-sorted set of path term classes, as the quotient of $\mathrm{PT}_{\boldsymbol{\mathcal{A}}}$ by the restriction of $\Theta^{[1]}$ to it, thus obtaining a subquotient of 
$\mathrm{T}_{\Sigma^{\boldsymbol{\mathcal{A}}}}(X)$. Next we show that $[\mathrm{PT}_{\boldsymbol{\mathcal{A}}}]$ is equipped with a structure of $S$-sorted partial 
$\Sigma^{\boldsymbol{\mathcal{A}}}$-algebra, a structure of $S$-sorted categorial $\Sigma$-algebra, and its coproduct, with a structure of Artinian order. We denote by 
$[\mathbf{PT}_{\boldsymbol{\mathcal{A}}}]$, $[\mathsf{PT}_{\boldsymbol{\mathcal{A}}}]$ and $(\coprod [\mathrm{PT}_{\boldsymbol{\mathcal{A}}}],\leq_{[\mathbf{PT}_{\boldsymbol{\mathcal{A}}}]})$, respectively, the corresponding $S$-sorted partial $\Sigma^{\boldsymbol{\mathcal{A}}}$-algebra, $S$-sorted categorial $\Sigma$-algebra and Artinian ordered set. Then, by means of the $S$-sorted mapping $\mathrm{CH}^{[1]}$ ---derived from $\mathrm{CH}^{(1)}$--- we prove one of the most important results of this part: That $[\mathbf{Pth}_{\boldsymbol{\mathcal{A}}}]$ and 
$[\mathbf{PT}_{\boldsymbol{\mathcal{A}}}]$ as well as $[\mathsf{Pth}_{\boldsymbol{\mathcal{A}}}]$ and $[\mathsf{PT}_{\boldsymbol{\mathcal{A}}}]$ and $(\coprod [\mathrm{Pth}_{\boldsymbol{\mathcal{A}}}],\leq_{[\mathbf{Pth}_{\boldsymbol{\mathcal{A}}}]})$ and $(\coprod [\mathrm{PT}_{\boldsymbol{\mathcal{A}}}],\leq_{[\mathbf{PT}_{\boldsymbol{\mathcal{A}}}]})$ are isomorphic. We also explicitly provide the inverse of $\mathrm{CH}^{[1]}$, which is the $S$-sorted mapping $\mathrm{ip}^{([1],X)@}$, derived from $\mathrm{ip}^{(1,X)@}$, the free completion of the $S$-sorted mapping $\mathrm{ip}^{(1,X)}$ from $X$ to $\mathrm{Pth}_{\boldsymbol{\mathcal{A}}}$. In this way we have established a Curry-Howard isomorphism type result in which the derivations, up to parallel reordering (of the classes of paths), admit an equivalent representation by means of terms, up to congruence (of the classes of path terms).

In addition, for a specification $\boldsymbol{\mathcal{E}}^{\boldsymbol{\mathcal{A}}}$ describing	 the properties that hold in the $S$-sorted partial $\Sigma^{\boldsymbol{\mathcal{A}}}$-algebra $[\mathbf{Pth}_{\boldsymbol{\mathcal{A}}}]$, and with regard to the $\mathrm{QE}$-variety determined by $\boldsymbol{\mathcal{E}}^{\boldsymbol{\mathcal{A}}}$, we prove that the $S$-sorted partial $\Sigma^{\boldsymbol{\mathcal{A}}}$-algebras   
$\mathbf{T}_{\boldsymbol{\mathcal{E}}^{\boldsymbol{\mathcal{A}}}}(\mathbf{Pth}_{\boldsymbol{\mathcal{A}}})$ and 
$[\mathbf{Pth}_{\boldsymbol{\mathcal{A}}}]$ are isomorphic, hence that $[\mathbf{PT}_{\boldsymbol{\mathcal{A}}}]$ and $\mathbf{T}_{\boldsymbol{\mathcal{E}}^{\boldsymbol{\mathcal{A}}}}(\mathbf{Pth}_{\boldsymbol{\mathcal{A}}})$ are also isomorphic, where $\mathbf{T}_{\boldsymbol{\mathcal{E}}^{\boldsymbol{\mathcal{A}}}}$ is the left adjoint of the canonical inclusion of $\mathsf{PAlg}(\boldsymbol{\mathcal{E}}^{\boldsymbol{\mathcal{A}}})$ into $\mathsf{PAlg}(\Sigma^{\boldsymbol{\mathcal{A}}})$. 

We finally define the notions of first-order elementary translation and of first-order translation and show that if two path terms $M$ and $N$ are such that the paths obtained from them under the action of $\mathrm{ip}^{(1,X)@}$ have the same $(0,1)$-source and $(0,1)$-target, then, for a first-order translation $T^{(1)}$, we have that $T^{(1)}(M)$ is a path term if and only if $T^{(1)}(N)$ is a path term and if either one of the two is a path term, then, when transformed into paths by means of $\mathrm{ip}^{(1,X)@}$, they still have the same $(0,1)$-source and $(0,1)$-target. Moreover, we show that if $M$ and $N$ are $\Theta^{[1]}$-related and either $T^{(1)}(M)$ or $T^{(1)}(N)$ is a path term, then $T^{(1)}(M)$ and $T^{(1)}(N)$ are also $\Theta^{[1]}$-related. These results are fundamental and constitute the basis for the definition of the notion of second-order path on path terms which will be the starting point of the second part.

%
%
%
%
%
%
%

\subsection{Second-order structures associated to a second-order rewriting system}
\textsf{Part~2} begins by defining, from a many-sorted term rewriting system $\boldsymbol{\mathcal{A}}$, a second-order many-sorted term rewriting system as an ordered pair $\boldsymbol{\mathcal{A}}^{(2)} = (\boldsymbol{\mathcal{A}},\mathcal{A}^{(2)})$, where $\mathcal{A}^{(2)}$ is a subset of $[\mathrm{PT}_{\boldsymbol{\mathcal{A}}}]^{2}$ subject to fulfilling a certain compatibility condition on sources and targets. Then we define the notion of second-order path in $\boldsymbol{\mathcal{A}}^{(2)}$ from a path term class to another. 

Before proceeding any further, it seems appropriate to point out that the passage from the first level, i.e., the world of many-sorted term rewriting systems, to the second level, i.e., the world of second-order many-sorted term rewriting systems, is similar to the passage from a continuous mapping between topological spaces to that of a homotopy between continuous mappings between topological spaces. And it is this topological analogy and its iteration (homotopies between homotopies, etc.), which will, to a large extent, guide us in our work. 

Next we show that second-order paths have a recursive decomposition, we call this the second-order path extraction algorithm. Following this, we define the $S$-sorted set $\mathrm{Pth}_{\boldsymbol{\mathcal{A}}^{(2)}}$, of second-order paths in $\boldsymbol{\mathcal{A}}^{(2)}$, and equip it with a structure of $\Sigma$-algebra. 
Moreover, we equip $\coprod\mathrm{Pth}_{\boldsymbol{\mathcal{A}}^{(2)}}$, the coproduct of $\mathrm{Pth}_{\boldsymbol{\mathcal{A}}^{(2)}}$, with a structure of Artinian order. After this, from $\boldsymbol{\mathcal{A}}^{(2)}$, we obtain an $S$-sorted signature $\Sigma^{\boldsymbol{\mathcal{A}}^{(2)}}$, which is an extension of the $S$-sorted signature $\Sigma^{\boldsymbol{\mathcal{A}}}$, and equip $\mathrm{Pth}_{\boldsymbol{\mathcal{A}}^{(2)}}$ with a structure of $S$-sorted partial $\Sigma^{\boldsymbol{\mathcal{A}}^{(2)}}$-algebra, which we denote by $\mathbf{Pth}_{\boldsymbol{\mathcal{A}}^{(2)}}$.

Afterwards, we define an $S$-sorted mapping $\mathrm{CH}^{(2)}$, the  $S$-sorted second-order Curry-Howard mapping, from the $S$-sorted set $\mathrm{Pth}_{\boldsymbol{\mathcal{A}}^{(2)}}$ to $\mathrm{T}_{\Sigma^{\boldsymbol{\mathcal{A}}^{(2)}}}(X)$, the underlying $S$-sorted set of the free $S$-sorted $\Sigma^{\boldsymbol{\mathcal{A}}^{(2)}}$-algebra $\mathbf{T}_{\Sigma^{\boldsymbol{\mathcal{A}}^{(2)}}}(X)$ on $X$. We  prove that the second-order Curry-Howard mapping is neither a
$\Sigma^{\boldsymbol{\mathcal{A}}^{(2)}}$-homomorphism, nor a
$\Sigma^{\boldsymbol{\mathcal{A}}}$-homomorphism from $\mathbf{Pth}_{\boldsymbol{\mathcal{A}}^{(2)}}$ to $\mathbf{T}_{\Sigma^{\boldsymbol{\mathcal{A}}^{(2)}}}(X)$.

We consider the quotient $[\mathrm{Pth}_{\boldsymbol{\mathcal{A}}^{(2)}}]$ of $\mathrm{Pth}_{\boldsymbol{\mathcal{A}}^{(2)}}$ by $\mathrm{Ker}(\mathrm{CH}^{(2)})$ and prove that it is equipped with a structure of $S$-sorted partial $\Sigma^{\boldsymbol{\mathcal{A}}^{(2)}}$-algebra, a structure of $S$-sorted categorial $\Sigma$-algebra and its coproduct, with a structure of Artinian order. We denote by $[\mathbf{Pth}_{\boldsymbol{\mathcal{A}}^{(2)}}]$, $[\mathsf{Pth}_{\boldsymbol{\mathcal{A}}^{(2)}}]$ and 
$(\coprod [\mathrm{Pth}_{\boldsymbol{\mathcal{A}}^{(2)}}],\leq_{[\mathbf{Pth}_{\boldsymbol{\mathcal{A}}^{(2)}}]})$, respectively, the corresponding partial $\Sigma^{\boldsymbol{\mathcal{A}}^{(2)}}$-algebra, $S$-sorted categorial $\Sigma$-algebra and Artinian ordered set. Following this, we show that $[\mathrm{Pth}_{\boldsymbol{\mathcal{A}}^{(2)}}]$ does not have a structure of $2$-category. To solve this problem, we define a new congruence $\Upsilon^{[1]}$ on $\mathbf{Pth}_{\boldsymbol{\mathcal{A}}^{(2)}}$ from which we get the partial 
$\Sigma^{\boldsymbol{\mathcal{A}}^{(2)}}$-algebra $[\mathbf{Pth}_{\boldsymbol{\mathcal{A}}^{(2)}}]_{\Upsilon^{[1]}} $, given by the quotient of $\mathbf{Pth}_{\boldsymbol{\mathcal{A}}^{(2)}}$ by $\Upsilon^{[1]}$. Following this we prove that $\coprod [\mathrm{Pth}_{\boldsymbol{\mathcal{A}}^{(2)}}]_{\Upsilon^{[1]}}$ is equipped with a structure of Artinian preorder, which we denote by $(\coprod [\mathrm{Pth}_{\boldsymbol{\mathcal{A}}^{(2)}}]_{\Upsilon^{[1]}},\leq_{[\mathbf{Pth}_{\boldsymbol{\mathcal{A}}^{(2)}}]_{\Upsilon^{[1]}}})$. 
We next consider $\llbracket \mathrm{Pth}_{\boldsymbol{\mathcal{A}}^{(2)}}\rrbracket$, the quotient of $\mathrm{Pth}_{\boldsymbol{\mathcal{A}}^{(2)}}$ by $\mathrm{Ker}(\mathrm{CH}^{(2)})\vee \Upsilon^{[1]}$, the supremum congruence of $\mathrm{Ker}(\mathrm{CH}^{(2)})$ and $\Upsilon^{[1]}$, and prove that it is equipped with a structure of $S$-sorted partial $\Sigma^{\boldsymbol{\mathcal{A}}^{(2)}}$-algebra,   a structure of $S$-sorted $2$-categorial $\Sigma$-algebra and, its coproduct, with a structure of Artinian preorder $\leq_{\llbracket \mathbf{Pth}_{\boldsymbol{\mathcal{A}}^{(2)}} \rrbracket}$. We denote by $\llbracket\mathbf{Pth}_{\boldsymbol{\mathcal{A}}^{(2)}}\rrbracket$, $\llbracket \mathsf{Pth}_{\boldsymbol{\mathcal{A}}^{(2)}} \rrbracket$ and $(\coprod\llbracket \mathrm{Pth}_{\boldsymbol{\mathcal{A}}^{(2)}} \rrbracket, \leq_{\llbracket \mathbf{Pth}_{\boldsymbol{\mathcal{A}}^{(2)}} \rrbracket})$, respectively, the corresponding partial $\Sigma^{\boldsymbol{\mathcal{A}}^{(2)}}$-algebra, $S$-sorted $2$-categorial $\Sigma$-algebra and
Artinian preordered set. This is another difference from the first part: what in the first part was an Artinian order is now only an Artinian preorder.

Afterwards, we define  a congruence $\Theta^{[2]}$ on $\mathbf{T}_{\Sigma^{\boldsymbol{\mathcal{A}}^{(2)}}}(X)$ and, from it, we obtain the $S$-sorted subset $\mathrm{PT}_{{\boldsymbol{\mathcal{A}}^{(2)}}}$, of second-order path terms, of $\mathbf{T}_{\Sigma^{\boldsymbol{\mathcal{A}}^{(2)}}}(X)$ as the $\Theta^{[2]}$-saturation of $\mathrm{CH}^{(2)}[\mathrm{Pth}_{\boldsymbol{\mathcal{A}}^{(2)}}]$. Then we define $[\mathrm{PT}_{{\boldsymbol{\mathcal{A}}^{(2)}}}]$, the $S$-sorted set of second-order path term classes, as the quotient of $\mathrm{PT}_{\boldsymbol{\mathcal{A}}^{(2)}}$ by the restriction of $\Theta^{[2]}$ to it, thus obtaining a subquotient of 
$\mathrm{T}_{\Sigma^{\boldsymbol{\mathcal{A}}^{(2)}}}(X)$. Next we show that $[\mathrm{PT}_{\boldsymbol{\mathcal{A}}^{(2)}}]$, is equipped with a structure of $S$-sorted partial $\Sigma^{\boldsymbol{\mathcal{A}}^{(2)}}$-algebra, a structure of $S$-sorted categorial $\Sigma$-algebra and, the coproduct of its underlying $S$-sorted set, with a structure of Artinian order. As it happened before with the quotient of second-order paths $[\mathrm{Pth}_{\boldsymbol{\mathcal{A}}^{(2)}}]$, we show that $[\mathrm{PT}_{\boldsymbol{\mathcal{A}}^{(2)}}]$ does not have a structure of $2$-category. To solve this problem, we define a new congruence $\Psi^{[1]}$ on $\mathbf{T}_{\boldsymbol{\mathcal{A}}^{(2)}}(X)$, from which we get the partial 
$\Sigma^{\boldsymbol{\mathcal{A}}^{(2)}}$-algebra $[\mathbf{PT}_{\boldsymbol{\mathcal{A}}^{(2)}}]_{\Psi^{[1]}} $, given by the quotient of $\mathbf{PT}_{\boldsymbol{\mathcal{A}}^{(2)}}$ by $\Psi^{[1]}$ restricted to it. Then we prove that $\coprod [\mathrm{PT}_{\boldsymbol{\mathcal{A}}^{(2)}}]_{\Psi^{[1]}}$ is equipped with a structure of Artinian preorder, that we denote by $(\coprod [\mathrm{PT}_{\boldsymbol{\mathcal{A}}^{(2)}}]_{\Psi^{[1]}},\leq_{[\mathbf{PT}_{\boldsymbol{\mathcal{A}}^{(2)}}]_{\Psi^{[1]}}})$. 
We next consider $\llbracket \mathrm{PT}_{\boldsymbol{\mathcal{A}}^{(2)}}\rrbracket$, the quotient of $\mathrm{PT}_{\boldsymbol{\mathcal{A}}^{(2)}}$ by $\Theta^{[2]}\vee \Psi^{[1]}$, the supremum of $\Theta^{[2]}$ and $\Psi^{[1]}$, simply denoted by $\Theta^{\llbracket 2 \rrbracket}$,  and prove that it is equipped with a structure of $S$-sorted partial $\Sigma^{\boldsymbol{\mathcal{A}}^{(2)}}$-algebra,   a structure of $S$-sorted $2$-categorial $\Sigma$-algebra and, its coproduct, with a structure of Artinian preorder $\leq_{\llbracket \mathbf{PT}_{\boldsymbol{\mathcal{A}}^{(2)}} \rrbracket}$. We denote by $\llbracket\mathbf{PT}_{\boldsymbol{\mathcal{A}}^{(2)}}\rrbracket$, $\llbracket \mathsf{PT}_{\boldsymbol{\mathcal{A}}^{(2)}} \rrbracket$ and $(\coprod\llbracket \mathrm{PT}_{\boldsymbol{\mathcal{A}}^{(2)}} \rrbracket, \leq_{\llbracket \mathbf{PT}_{\boldsymbol{\mathcal{A}}^{(2)}} \rrbracket})$, respectively, the corresponding partial $\Sigma^{\boldsymbol{\mathcal{A}}^{(2)}}$-algebra, $S$-sorted $2$-categorial $\Sigma$-algebra and
Artinian preordered set. 

Then, by means of the $S$-sorted mapping $\mathrm{CH}^{\llbracket 2\rrbracket}$ ---derived from $\mathrm{CH}^{(2)}$--- we prove one of the most important results of the second part: That $\llbracket \mathbf{Pth}_{\boldsymbol{\mathcal{A}}^{(2)}}\rrbracket$ and 
$\llbracket \mathbf{PT}_{\boldsymbol{\mathcal{A}}^{(2)}}\rrbracket$ as well as $\llbracket \mathsf{Pth}_{\boldsymbol{\mathcal{A}}^{(2)}}\rrbracket$ and 
$\llbracket \mathsf{PT}_{\boldsymbol{\mathcal{A}}^{(2)}}\rrbracket$ and $(\coprod \llbracket \mathbf{Pth}_{\boldsymbol{\mathcal{A}}^{(2)}}\rrbracket,\leq_{\llbracket \mathbf{Pth}_{\boldsymbol{\mathcal{A}}^{(2)}}\rrbracket})$ and $(\coprod \llbracket \mathbf{PT}_{\boldsymbol{\mathcal{A}}^{(2)}}\rrbracket,\leq_{\llbracket \mathbf{PT}_{\boldsymbol{\mathcal{A}}^{(2)}}\rrbracket})$ are isomorphic. We also explicitly provide the inverse of $\mathrm{CH}^{\llbracket 2 \rrbracket}$, which is the $S$-sorted mapping $\mathrm{ip}^{(\llbracket 2\rrbracket,X)@}$, derived from $\mathrm{ip}^{(2,X)@}$, the free completion of the $S$-sorted mapping $\mathrm{ip}^{(2,X)}$ from $X$ to $\mathrm{Pth}_{\boldsymbol{\mathcal{A}}^{(2)}}$. In this way we have established a second-order Curry-Howard type result in which the second-order derivations, up to parallel reordering and second-order categorial concerns (of the classes of second-order paths), admit an equivalent representation by means of terms, up to congruence (of the classes of second-order path terms).

\begin{figure}
\begin{center}
\begin{adjustbox}{width=\textwidth}
\begin{tikzpicture}
[ACliment/.style={-{To [angle'=45, length=5.75pt, width=4pt, round]}}, 
scale=.8
]

\node[] (T) at (0,0) {$\mathbf{T}_{\Sigma}(X)$};

\node[] (Pth) at (4,0) {$\mathrm{Pth}_{\boldsymbol{\mathcal{A}}}$};
\node[] (PT) at (6,2) {$\mathrm{PT}_{\boldsymbol{\mathcal{A}}}$};
\node[] (PthQ) at (6,-2) {$[\mathrm{Pth}_{\boldsymbol{\mathcal{A}}}]$};
\node[] (PTQ) at (8,0) {$[\mathrm{PT}_{\boldsymbol{\mathcal{A}}}]$};
\node[] () at (7,-1) {\rotatebox{45}{${\scriptsize\cong}$}};

\node[] (Pth2) at (13,0) {$\mathrm{Pth}_{\boldsymbol{\mathcal{A}}^{(2)}}$};
\node[] (PT2) at (15,2) {$\mathrm{PT}_{\boldsymbol{\mathcal{A}}^{(2)}}$};
\node[] (PthQ2) at (15,-2) {$\llbracket\mathrm{Pth}_{\boldsymbol{\mathcal{A}}^{(2)}}\rrbracket$};
\node[] (PTQ2) at (17,0) {$\llbracket\mathrm{PT}_{\boldsymbol{\mathcal{A}}^{(2)}}\rrbracket$};
\node[] () at (16,-1) {\rotatebox{45}{${\scriptsize\cong}$}};

\draw[ACliment]  ($(T)+(.7,0)$) to node [midway, fill=white] 
{$\scriptsize\mathrm{ip}^{(1,0)\sharp}$} ($(Pth)+(-.5,0)$);
\draw[ACliment]  ($(Pth)+(-.5,.3)$) to node [above, pos=.5] 
{$\scriptsize\mathrm{sc}^{(0,1)}$} ($(T)+(.7,.3)$);
\draw[ACliment]  ($(Pth)+(-.5,-.3)$) to node [below, pos=.5] 
{$\scriptsize\mathrm{tg}^{(0,1)}$} ($(T)+(.7,-.3)$);

\draw[ACliment]  ($(PTQ)+(.6,0)$) to node [midway, fill=white] 
{$\scriptsize\mathrm{ip}^{(\llbracket 2\rrbracket,[1])\sharp}$} ($(Pth2)+(-.8,0)$);
\draw[ACliment]  ($(Pth2)+(-.8,.3)$) to node [above, pos=.5] 
{$\scriptsize\mathrm{sc}^{([1],\llbracket 2\rrbracket)}$} ($(PTQ)+(.6,.3)$);
\draw[ACliment]  ($(Pth2)+(-.8,-.3)$) to node [below, pos=.5] 
{$\scriptsize\mathrm{tg}^{([1],\llbracket 2\rrbracket)}$} ($(PTQ)+(.6,-.3)$);

\draw[ACliment]  ($(Pth)+(.1,.3)$) to node [left, pos=.5] 
{$\scriptsize\mathrm{CH}^{(1)}$} ($(PT)+(-.4,-.2)$);
\draw[ACliment]  ($(PT)+(-.2,-.4)$) to node [right, pos=.5] 
{$\scriptsize\mathrm{ip}^{(1,X)@}$} ($(Pth)+(.3,.1)$);

\draw[ACliment]  ($(PthQ)+(.1,.3)$) to node [left, pos=.5] 
{$\scriptsize\mathrm{CH}^{[1]}$} ($(PTQ)+(-.4,-.2)$);
\draw[ACliment]  ($(PTQ)+(-.2,-.4)$) to node [right, pos=.5] 
{$\scriptsize\mathrm{ip}^{([1],X)@}$} ($(PthQ)+(.3,.1)$);

\draw[ACliment]  (Pth) to node [below left, pos=.5] 
{$\scriptsize\mathrm{pr}^{\mathrm{Ker}(\mathrm{CH}^{(1)})}$} (PthQ);
\draw[ACliment]  (PT) to node [above right, pos=.5] 
{$\scriptsize\mathrm{pr}^{\Theta^{[1]}}$} (PTQ);

\draw[ACliment]  ($(Pth2)+(.1,.3)$) to node [left, pos=.5] 
{$\scriptsize\mathrm{CH}^{(2)}$} ($(PT2)+(-.4,-.2)$);
\draw[ACliment]  ($(PT2)+(-.2,-.4)$) to node [right, pos=.5] 
{$\scriptsize\mathrm{ip}^{(2,X)@}$} ($(Pth2)+(.3,.1)$);

\draw[ACliment]  ($(PthQ2)+(.1,.3)$) to node [left, pos=.5] 
{$\scriptsize\mathrm{CH}^{\llbracket 2\rrbracket}$} ($(PTQ2)+(-.4,-.2)$);
\draw[ACliment]  ($(PTQ2)+(-.2,-.4)$) to node [right, pos=.5] 
{$\scriptsize\mathrm{ip}^{(\llbracket 2\rrbracket,X)@}$} ($(PthQ2)+(.3,.1)$);

\draw[ACliment]  (Pth2) to node [below left, pos=.5] 
{$\scriptsize\mathrm{pr}^{\mathrm{Ker}(\mathrm{CH}^{(2)})}$} (PthQ2);
\draw[ACliment]  (PT2) to node [above right, pos=.5] 
{$\scriptsize\mathrm{pr}^{\Theta^{\llbracket 2 \rrbracket}}$} (PTQ2);

\end{tikzpicture}
\end{adjustbox}
\caption{Representation of the first two parts and their interrelationships.}
\label{FA0}
\end{center}
\end{figure}


In addition, for a specification $\boldsymbol{\mathcal{E}}^{\boldsymbol{\mathcal{A}}^{(2)}}$ describing the properties that hold in the $S$-sorted partial $\Sigma^{\boldsymbol{\mathcal{A}}^{(2)}}$-algebra $\llbracket \mathbf{Pth}_{\boldsymbol{\mathcal{A}}^{(2)}}\rrbracket$, and with regard to the $\mathrm{QE}$-variety determined by $\boldsymbol{\mathcal{E}}^{\boldsymbol{\mathcal{A}}^{(2)}}$, we prove that the $S$-sorted partial $\Sigma^{\boldsymbol{\mathcal{A}}^{(2)}}$-algebras   
$\mathbf{T}_{\boldsymbol{\mathcal{E}}^{\boldsymbol{\mathcal{A}}^{(2)}}}(\mathbf{Pth}_{\boldsymbol{\mathcal{A}}^{(2)}})$ and 
$\llbracket\mathbf{Pth}_{\boldsymbol{\mathcal{A}}^{(2)}}\rrbracket$ are isomorphic, hence that $\llbracket\mathbf{PT}_{\boldsymbol{\mathcal{A}}^{(2)}}\rrbracket$ and $\mathbf{T}_{\boldsymbol{\mathcal{E}}^{\boldsymbol{\mathcal{A}}^{(2)}}}(\mathbf{Pth}_{\boldsymbol{\mathcal{A}}^{(2)}})$ are also isomorphic, where $\mathbf{T}_{\boldsymbol{\mathcal{E}}^{\boldsymbol{\mathcal{A}}^{(2)}}}$ is the left adjoint of the canonical inclusion of $\mathsf{PAlg}(\boldsymbol{\mathcal{E}}^{\boldsymbol{\mathcal{A}}^{(2)}})$ into $\mathsf{PAlg}(\Sigma^{\boldsymbol{\mathcal{A}}^{(2)}})$.

We finally define the notions of second-order elementary translation and of second-order translation and show that if two second-order path terms $M$ and $N$ are such that the second-order paths obtained from them under the action of $\mathrm{ip}^{(2,X)@}$ have the same $([1],2)$-source and $([1],2)$-target, then, for a second-order translation $T^{(2)}$, we have that $T^{(2)}(M)$ is a second-order path term if and only if $T^{(2)}(N)$ is a second-order path term and if either one of the two is a second-order path term, then, when transformed into second-order paths by means of $\mathrm{ip}^{(2,X)@}$, they still have the same $([1],2)$-source and $([1],2)$-target. Moreover, we show that if $M$ and $N$ are $\Theta^{\llbracket 2 \rrbracket}$-related and either $T^{(2)}(M)$ or $T^{(2)}(N)$ is a second-order path term, then $T^{(2)}(M)$ and $T^{(2)}(N)$ are also $\Theta^{\llbracket 2 \rrbracket}$-related. These results are foundational and represent the extension of the aforementioned results to higher-order instances. 

Let us point out that it is at the second level where the structural richness, that was hidden at the first level due to the one-dimensional constraint, unfolds and the pattern for recursively defining a tower associated to a higher-order many-sorted term rewriting system can, at last, be discerned.

\subsection{Morphisms of rewriting systems}

In \textsf{Part~3}, we start by defining $\mathsf{Rws}_{\mathfrak{d}}^{(0)}$,  the category of zeroth order many-sorted rewriting systems, which is at the basis of all the work. Let $\mathsf{Sig}_{\boldsymbol{\mathfrak{d}}}$ be the category whose objects are the many-sorted signatures, i.e., the ordered pairs $(S,\Sigma)$, where $S$ is a set, the set of sorts, and $\Sigma$ an $S$-sorted signature, i.e., a mapping from $S^{\star}\times S$ to a Grothendieck universe $\boldsymbol{\mathcal{U}}$, where $S^{\star}$ is the underlying set of the free monoid on $S$, and whose morphisms from the many-sorted signature $(S,\Sigma)$ to  the many-sorted signature $(T,\Lambda)$ are the derivors $\mathbf{d} = (\varphi,d)$ (see~\cite{CVCL18}). 
Let $\mathsf{Alg}_{\mathfrak{d}}$ be the category whose objects are the many-sorted algebras, i.e., the ordered pairs $((S,\Sigma),\mathbf{A})$, where $(S,\Sigma)$ is a many-sorted signature and $\mathbf{A}$ a $\Sigma$-algebra, and whose morphisms from the many-sorted algebra $((S,\Sigma),\mathbf{A})$ to the many-sorted algebra $((T,\Lambda),\mathbf{B})$ are the ordered pairs $(\mathbf{d},f)$, with $\mathbf{d} = (\varphi,d)$ a derivor from $(S,\Sigma)$ to $(T,\Lambda)$ and $f$ a homomorphism of $\Sigma$-algebras from $\mathbf{A}$ to $\mathbf{d}^{\ast}_{\mathfrak{d}}(\mathbf{B})$, where $\mathbf{d}^{\ast}_{\mathfrak{d}}$ is the functor from $\mathsf{Alg}(\Lambda)$ to $\mathsf{Alg}(\Sigma)$ determined by the derivor $\mathbf{d}$ (see~\cite{CVCL18}). We remark that from $\mathsf{Sig}_{\mathfrak{d}}$ to $\mathsf{Cat}$, the category of $\boldsymbol{\mathcal{U}}$-categories, there exists a contravariant functor, denoted by $\mathrm{Alg}_{\mathfrak{d}}$, and that the category $\mathsf{Alg}_{\mathfrak{d}}$ is $\int^{\mathsf{Sig}_{\scriptscriptstyle\mathfrak{d}}}\mathrm{Alg}_{\mathfrak{d}}$, i.e., the category obtained by means of the Grothendieck construction applied to $\mathrm{Alg}_{\mathfrak{d}}$.

Then $\mathsf{Rws}_{\mathfrak{d}}^{(0)}$ is the category whose objects are the ordered triples $(S,\Sigma,X)$, where $(S,\Sigma)$ is a many-sorted signature and $X$ an $S$-sorted set, which we will call zeroth order many-sorted rewriting systems, 
and whose morphisms from the zeroth order many-sorted rewriting system $\boldsymbol{\mathcal{A}}^{(0)}=(S,\Sigma,X)$ to the zeroth order many-sorted rewriting system $\boldsymbol{\mathcal{B}}^{(0)}=(T,\Lambda,Y)$ are the ordered triples $((S,\Sigma,X),\mathbf{f}^{(0)},(T,\Lambda,Y))$, denoted by $\mathbf{f}^{(0)}\colon \boldsymbol{\mathcal{A}}^{(0)}\mor \boldsymbol{\mathcal{B}}^{(0)}$ for short, in which $\mathbf{f}^{(0)}=(\varphi,d,f^{(0)})$, where 
\begin{enumerate}
\item $(\varphi,d)$ is a derivor from $(S,\Sigma)$ to $(T,\Lambda)$ and 
\item $f^{(0)}$ is a mapping from $X$ to $\mathrm{T}_{\Lambda}(Y)_{\varphi}=(\mathrm{T}_{\Lambda}(Y)_{\varphi(s)})_{s\in S}$.
\end{enumerate}

Note that $\mathrm{T}_{\Lambda}(Y)_{\varphi}$ is the underlying $S$-sorted set of the $\Sigma$-algebra $\mathbf{d}^{\ast}_{\mathfrak{d}}(\mathbf{T}_{\Lambda}(Y))$. Let us also note that there exists a full embedding of $\mathsf{Rws}_{\mathfrak{d}}^{(0)}$ into $\mathsf{Alg}_{\mathfrak{d}}$, because to give an $S$-sorted mapping $f^{(0)}$ from $X$ to $\mathrm{T}_{\Lambda}(Y)_{\varphi}$ is, naturally, equivalent to give a $\Sigma$-homomorphism $f^{(0)\sharp}$ from $\mathbf{T}_{\Sigma}(X)$ to $\mathbf{d}^{\ast}_{\mathfrak{d}}(\mathbf{T}_{\Lambda}(Y))$, obtained by the universal property of $\mathrm{T}_{\Sigma}(X)$. Actually, we prove that $\mathsf{Rws}_{\mathfrak{d}}^{(0)}$ is isomorphic to $\mathsf{Tw}^{(0)}_{\mathfrak{d}}$, the full subcategory of 
$\mathsf{Alg}_{\mathfrak{d}}$ determined by the many-sorted algebras of the form $(\mathbcal{A}^{(0)},\mathbf{T}_{\Sigma}(X))$, with $\mathbcal{A}^{(0)}$ varying over the zeroth-order many-sorted rewriting systems, what we have called in this work zeroth-order towers. 

We next consider first-order rewriting systems and first-order morphisms between them. Let $\boldsymbol{\mathcal{A}}^{(1)} = (S,\Sigma,X,\mathcal{A}^{(1)})=(\boldsymbol{\mathcal{A}}^{(0)},\mathcal{A}^{(1)})$ and $\boldsymbol{\mathcal{B}}^{(1)} = (T,\Lambda,Y,\mathcal{B}^{(1)})=(\boldsymbol{\mathcal{B}}^{(0)},\mathcal{B}^{(1)})$ be first-order many-sorted rewriting systems. A morphism from $\boldsymbol{\mathcal{A}}^{(1)}$ to $\boldsymbol{\mathcal{B}}^{(1)}$ will be an ordered triple 
$(\mathbcal{A}^{(1)},\mathbf{f}^{(1)},\mathbcal{B}^{(1)})$, denoted by 
$\mathbf{f}^{(1)}\colon \boldsymbol{\mathcal{A}}^{(1)}\mor \boldsymbol{\mathcal{B}}^{(1)}$ or $\mathbf{f}^{(1)}$ 
for short, in which $\mathbf{f}^{(1)}$ is an ordered pair $(\mathbf{f}^{(0)},f^{(1)})$ where 
\begin{enumerate}
\item $\mathbf{f}^{(0)} = (\varphi,d,f^{(0)})$ is a morphism from $\boldsymbol{\mathcal{A}}^{(0)}$ to 
$\boldsymbol{\mathcal{B}}^{(0)}$ and
\item $f^{(1)}$ is a mapping from $\mathcal{A}^{(1)}$ to 
$\mathrm{Pth}_{\mathbcal{B}^{(1)},\varphi} = (\mathrm{Pth}_{\mathbcal{B}^{(1)},\varphi(s)})_{s\in S}$ such that, for every $s\in S$ and every $\mathfrak{p} = (M,N)\in \mathcal{A}^{(1)}_{s}$, we have that 
$$
f^{(1)}_{s}(\mathfrak{p})\in \mathrm{Pth}_{\mathbcal{B}^{(1)},\varphi(s)}(f^{(0)\sharp}_{s}(M),f^{(0)\sharp}_{s}(N)).
$$
\end{enumerate}

From this definition we obtain, by recursion on $(\coprod \mathrm{Pth}_{\mathbcal{A}^{(1)}},\leq_{\mathbf{Pth}_{\mathbcal{A}^{(1)}}})$, a unique $\Sigma$-homomor\-phism $f^{(1)\flat}$ from $\mathbf{Pth}_{\mathbcal{A}^{(1)}}^{(0,1)}$ to $\mathbf{d}^{\ast}_{\mathfrak{d}}(\mathbf{Pth}_{\mathbcal{B}^{(1)}}^{(0,1)})$ such that $f^{(1)\flat}$ coincides with $f^{(0)\sharp}$ when restricted to the $(1,0)$-identity paths and with $f^{(1)}$ when restricted to the echelons, i.e., to the one-step paths canonically associated to the rewrite rules on $\mathcal{A}^{(1)}$. We next provide $[\mathrm{Pth}_{\boldsymbol{\mathcal{B}}^{(1)}}]_{\varphi}=([\mathrm{Pth}_{\boldsymbol{\mathcal{B}}^{(1)}}]_{\varphi(s)})_{s\in S}$, i.e., the underlying $S$-sorted set of the $\Sigma$-algebra  $\mathbf{d}^{\ast}_{\mathfrak{d}}([\mathbf{Pth}_{\mathbcal{B}^{(1)}}^{(0,1)}])$, with a structure of partial $\Sigma^{\boldsymbol{\mathcal{A}}^{(1)}}$-algebra, that we denote by $[\mathbf{Pth}^{\mathbf{f}^{(1)}}_{\boldsymbol{\mathcal{B}}^{(1)}}]$, and we prove that it belongs to the QE-variety determined by $\mathbcal{E}^{\mathbcal{A}^{(1)}}$. Since the partial $\Lambda^{\mathbcal{B}^{(1)}}$-algebras $[\mathbf{PT}_{\boldsymbol{\mathcal{B}}^{(1)}}]$ and $\mathbf{T}_{\mathbcal{E}^{\boldsymbol{\mathcal{B}}^{(1)}}}(\mathbf{Pth}_{\mathbcal{B}^{(1)}})$ are isomorphic to  $[\mathbf{Pth}_{\boldsymbol{\mathcal{B}}^{(1)}}]$, an analogous statement holds for the respective partial $\Sigma^{\boldsymbol{\mathcal{A}}^{(1)}}$-algebras $[\mathbf{PT}^{\mathbf{f}^{(1)}}_{\boldsymbol{\mathcal{B}}^{(1)}}]$ and $\mathbf{T}_{\mathbcal{E}^{\boldsymbol{\mathcal{B}}^{(1)}}}^{\mathbf{f}^{(1)}}(\mathbf{Pth}_{\mathbcal{B}^{(1)}})$, built following a similar approach.

Thus, since we have proved that $[\mathbf{Pth}_{\boldsymbol{\mathcal{A}}^{(1)}}]$ is isomorphic to $\mathbf{T}_{\mathbcal{E}^{\boldsymbol{\mathcal{A}}^{(1)}}}(\mathbf{Pth}_{\mathbcal{A}})$ (or to $[\mathbf{PT}_{\mathbcal{A}^{(1)}}]$), we obtain, by the universal property of it, a unique $\Sigma^{\mathbcal{A}^{(1)}}$-homomorphism from $[\mathbf{Pth}_{\boldsymbol{\mathcal{A}}^{(1)}}]$ to $[\mathbf{Pth}_{\boldsymbol{\mathcal{B}}^{(1)}}^{\mathbf{f}^{(1)}}]_{\varphi}$ (from $\mathbf{T}_{\mathbcal{E}^{\boldsymbol{\mathcal{A}}^{(1)}}}(\mathbf{Pth}_{\mathbcal{A}^{(1)}})$ to $\mathbf{T}_{\mathbcal{E}^{\boldsymbol{\mathcal{B}}^{(1)}}}^{\mathbf{f}^{(1)}}(\mathbf{Pth}_{\mathbcal{B}^{(1)}})$ or from $[\mathbf{PT}_{\mathbcal{A}^{(1)}}]$ to $[\mathbf{PT}_{\mathbcal{B}^{(1)}}^{\mathbf{f}^{(1)}}]$), which we denote by $f^{[1]@}$, satisfying $$f^{[1]@}\circ \mathrm{pr}^{\mathrm{Ker}(\mathrm{CH}^{(1)})}_{\mathbcal{A}^{(1)}}=\mathrm{pr}^{\mathrm{Ker}(\mathrm{CH}^{(1)})}_{\mathbcal{B}^{(1)}, \varphi}\circ f^{(1)\flat}.$$ 

Now, with respect to first-order many-sorted rewriting systems and morphisms on it, we prove that it does not form a category since, among other things, the natural composition of morphisms is not associative. To overcome this problem, we introduce an equivalence relation on the morphisms. Given two first-order morphisms $\mathbf{f}^{(1)}=(\mathbf{f}^{(0)},f^{(1)})$ and $\mathbf{g}^{(1)}=(\mathbf{g}^{(0)},g^{(1)})$, from $\mathbcal{A}^{(1)}$ to $\mathbcal{B}^{(1)}$, we say that $\mathbf{f}^{(1)}$ is equivalent to $\mathbf{g}^{(1)}$, written $\mathbf{f}^{(1)}\cong^{(1)} \mathbf{g}^{(1)}$, if and only if, $\mathbf{f}^{(0)}=\mathbf{g}^{(0)}$ and $\mathrm{pr}^{\mathrm{Ker}(\mathrm{CH}^{(1)})}_{\mathbcal{B}^{(1)},\varphi}\circ f^{(1)\flat}=\mathrm{pr}^{\mathrm{Ker}(\mathrm{CH}^{(1)})}_{\mathbcal{B}^{(1)},\varphi}\circ g^{(1)\flat}$. We prove that, for equivalent morphisms $\mathbf{f}^{(1)}\cong^{(1)} \mathbf{g}^{(1)}$, it holds that $f^{[1]@}=g^{[1]@}$ and the respective partial $\Sigma^{\mathbcal{A}^{(1)}}$-algebras $[\mathbf{Pth}^{\mathbf{f}^{(1)}}_{\mathbcal{B}^{(1)}}]$ and $[\mathbf{Pth}^{\mathbf{g}^{(1)}}_{\mathbcal{B}^{(1)}}]$ are isomorphic. This allows us to construct the category $\mathsf{Rws}_{\mathfrak{d}}^{[1]}$, whose objects are first-order many-sorted rewriting systems and morphisms are $\cong^{(1)}$-classes of morphisms of first-order many-sorted rewriting systems. Actually, we prove that $\mathsf{Rws}^{[1]}_{\mathfrak{d}}$ is a category and is isomorphic to the category of first-order towers, denoted by $\mathsf{Tw}^{[1]}_{\mathfrak{d}}$, determined by the many-sorted partial algebras of the form $(\mathbcal{A}^{(1)}, \mathbf{T}_{\mathbcal{E}^{\boldsymbol{\mathcal{A}}^{(1)}}}(\mathbf{Pth}_{\mathbcal{A}^{(1)}}))$, with $\mathbcal{A}^{(1)}$ varying over the first-order many-sorted rewriting systems.

Finally, we next consider second-order rewriting systems and second-order morphisms between them. Let $\boldsymbol{\mathcal{A}}^{(2)} = (S,\Sigma,X,\mathcal{A}^{(1)},\mathcal{A}^{(2)})=(\boldsymbol{\mathcal{A}}^{(1)},\mathcal{A}^{(2)})$ and $\boldsymbol{\mathcal{B}}^{(2)} = (T,\Lambda,Y,\mathcal{B}^{(1)},\mathcal{B}^{(2)})=(\boldsymbol{\mathcal{B}}^{(1)},\mathcal{B}^{(2)})$ be second-order many-sorted rewriting systems. A morphism from $\boldsymbol{\mathcal{A}}^{(2)}$ to $\boldsymbol{\mathcal{B}}^{(2)}$ will be an ordered triple 
$(\mathbcal{A}^{(2)},\mathbf{f}^{(2)},\mathbcal{B}^{(2)})$, denoted by 
$\mathbf{f}^{(2)}\colon \boldsymbol{\mathcal{A}}^{(2)}\mor \boldsymbol{\mathcal{B}}^{(2)}$ or $\mathbf{f}^{(2)}$ 
for short, in which $\mathbf{f}^{(2)}$ is an ordered pair $(\mathbf{f}^{(1)},f^{(2)})$ where 
\begin{enumerate}
\item $\mathbf{f}^{(1)} = (\mathbf{f}^{(0)},f^{(1)})$ is a morphism from $\boldsymbol{\mathcal{A}}^{(1)}$ to 
$\boldsymbol{\mathcal{B}}^{(1)}$ and
\item $f^{(2)}$ is a mapping from $\mathcal{A}^{(2)}$ to 
$\mathrm{Pth}_{\mathbcal{B}^{(2)},\varphi} = (\mathrm{Pth}_{\mathbcal{B}^{(2)},\varphi(s)})_{s\in S}$ such that, for every $s\in S$ and every $\mathfrak{p}^{(2)} = ([M]_{s},[N]_{s})\in \mathcal{A}^{(2)}_{s}$, we have that 
$$
f^{(2)}_{s}(\mathfrak{p}^{(2)})\in \mathrm{Pth}_{\mathbcal{B}^{(2)},\varphi(s)}(f^{[1]@}_{s}([M]_{s}),f^{[1]@}_{s}([N]_{s})).
$$
\end{enumerate}

From this definition we obtain, by recursion on $(\coprod \mathrm{Pth}_{\mathbcal{A}^{(2)}},\leq_{\mathbf{Pth}_{\mathbcal{A}^{(2)}}})$, a $\Sigma$-homomor\-phism $f^{(2)\flat}$ from $\mathbf{Pth}_{\mathbcal{A}^{(2)}}^{(0,2)}$ to $\mathbf{d}^{\ast}_{\mathfrak{d}}(\mathbf{Pth}_{\mathbcal{B}^{(2)}}^{(0,2)})$ such that $f^{(2)\flat}$ coincides with $f^{[1]@}$ when restricted to the $(2,[1])$-identity paths and with $f^{(2)}$ when restricted to the second-order echelons, i.e., to the one-step second-order paths canonically associated to the second-order rewrite rules on $\mathcal{A}^{(2)}$. We next provide $\llbracket\mathrm{Pth}_{\boldsymbol{\mathcal{B}}^{(2)}}\rrbracket_{\varphi}=(\llbracket\mathrm{Pth}_{\boldsymbol{\mathcal{B}}^{(2)}}\rrbracket_{\varphi(s)})_{s\in S}$, i.e., the underlying $S$-sorted set of the $\Sigma$-algebra  $\mathbf{d}^{\ast}_{\mathfrak{d}}(\llbracket\mathbf{Pth}_{\mathbcal{B}^{(2)}}^{(0,2)}\rrbracket)$, with a structure of partial $\Sigma^{\boldsymbol{\mathcal{A}}^{(2)}}$-algebra, that we denote by $\llbracket \mathbf{Pth}^{\mathbf{f}^{(2)}}_{\boldsymbol{\mathcal{B}}^{(2)}}\rrbracket$, and we prove that it belongs to the QE-variety determined by $\mathbcal{E}^{\mathbcal{A}^{(2)}}$. Since the partial $\Lambda^{\mathbcal{B}^{(2)}}$-algebras $\llbracket\mathbf{PT}_{\boldsymbol{\mathcal{B}}^{(2)}}\rrbracket$ and $\mathbf{T}_{\mathbcal{E}^{\boldsymbol{\mathcal{B}}^{(2)}}}(\mathbf{Pth}_{\mathbcal{B}^{(2)}})$ are isomorphic to  $\llbracket\mathbf{Pth}_{\boldsymbol{\mathcal{B}}^{(2)}}\rrbracket$, an analogous statement holds for the respective partial $\Sigma^{\boldsymbol{\mathcal{A}}^{(2)}}$-algebras $\llbracket\mathbf{PT}^{\mathbf{f}^{(2)}}_{\boldsymbol{\mathcal{B}}^{(2)}}\rrbracket$ and $\mathbf{T}_{\mathbcal{E}^{\boldsymbol{\mathcal{B}}^{(2)}}}^{\mathbf{f}^{(2)}}(\mathbf{Pth}_{\mathbcal{B}^{(2)}})$, built following a similar approach. 

Thus, since we have proved that $\llbracket\mathbf{Pth}_{\boldsymbol{\mathcal{A}}^{(2)}}\rrbracket$ is isomorphic to $\mathbf{T}_{\mathbcal{E}^{\boldsymbol{\mathcal{A}}^{(2)}}}(\mathbf{Pth}_{\mathbcal{A}^{(2)}})$ (or to $\llbracket\mathbf{PT}_{\mathbcal{A}^{(2)}}\rrbracket$), we obtain, by the universal property of it, a unique $\Sigma^{\mathbcal{A}^{(2)}}$-homomorphism from $\llbracket\mathbf{Pth}_{\boldsymbol{\mathcal{A}}^{(2)}}\rrbracket$ to $\llbracket\mathbf{Pth}_{\boldsymbol{\mathcal{B}}^{(2)}}^{\mathbf{f}^{(2)}}\rrbracket_{\varphi}$ (from $\mathbf{T}_{\mathbcal{E}^{\boldsymbol{\mathcal{A}}^{(2)}}}(\mathbf{Pth}_{\mathbcal{A}^{(2)}})$ to $\mathbf{T}_{\mathbcal{E}^{\boldsymbol{\mathcal{B}}^{(2)}}}^{\mathbf{f}^{(2)}}(\mathbf{Pth}_{\mathbcal{B}^{(2)}})$ or from $\llbracket \mathrm{PT}_{\mathbcal{A}^{(2)}}\rrbracket$ to $\llbracket \mathrm{PT}_{\mathbcal{B}^{(2)}}^{\mathbf{f}^{(2)}}\rrbracket$), which we denote by $f^{\llbracket 2 \rrbracket@}$, satisfying $$f^{\llbracket 2 \rrbracket@}\circ \mathrm{pr}^{\llbracket \cdot \rrbracket}_{\mathbcal{A}^{(2)}}=\mathrm{pr}^{\llbracket \cdot \rrbracket}_{\mathbcal{B}^{(2)}, \varphi}\circ f^{(2)\flat}.$$ 

Now, with respect to second-order many-sorted rewriting systems and morphisms on it, we prove that it does not form a category since, among other things, the natural composition of morphisms is not associative. To overcome this problem, we introduce an equivalence relation on the morphisms. Given two morphisms $\mathbf{f}^{(2)} = (\mathbf{f}^{(1)}, f^{(2)})$ and $\mathbf{g}^{(2)} = (\mathbf{g}^{(1)}, g^{(2)})$, from $\mathbcal{A}^{(2)}$ to $\mathbcal{B}^{(2)}$, we say that $\mathbf{f}^{(2)}$ is equivalent to $\mathbf{g}^{(2)}$, written $\mathbf{f}^{(2)}\cong^{(2)} \mathbf{g}^{(2)}$, if and only if, $\mathbf{f}^{(1)} \cong^{(1)} \mathbf{g}^{(1)}$ and $\mathrm{pr}^{\llbracket \cdot \rrbracket}_{\mathbcal{B}^{(2)},\varphi}\circ f^{(1)\flat}=\mathrm{pr}^{\llbracket \cdot \rrbracket}_{\mathbcal{B}^{(2)},\varphi}\circ g^{(2)\flat}$. We prove that, for equivalent morphisms $\mathbf{f}^{(2)}\cong^{(2)} \mathbf{g}^{(2)}$, it holds that $f^{\llbracket 2 \rrbracket @}=g^{\llbracket 2 \rrbracket@}$ and the respective partial $\Sigma^{\mathbcal{A}^{(2)}}$-algebras $\llbracket\mathbf{Pth}^{\mathbf{f}^{(2)}}_{\mathbcal{B}^{(2)}}\rrbracket$ and $\llbracket\mathbf{Pth}^{\mathbf{g}^{(2)}}_{\mathbcal{B}^{(2)}}\rrbracket$ are isomorphic. This allows us to construct the category $\mathsf{Rws}_{\mathfrak{d}}^{\llbracket 2 \rrbracket}$, whose objects are second-order many-sorted rewriting systems and morphisms are $\cong^{(2)}$-classes of morphisms of second-order many-sorted rewriting systems. Actually, we prove that $\mathsf{Rws}^{\llbracket 2 \rrbracket}_{\mathfrak{d}}$ is a category and is isomorphic to the category of second-order towers, denoted by $\mathsf{Tw}^{\llbracket 2 \rrbracket}_{\mathfrak{d}}$, determined by the many-sorted partial algebras of the form $(\mathbcal{A}^{(2)}, \mathbf{T}_{\mathbcal{E}^{\boldsymbol{\mathcal{A}}^{(2)}}}(\mathbf{Pth}_{\mathbcal{A}^{(2)}}))$, with $\mathbcal{A}^{(2)}$ varying over the second-order many-sorted rewriting systems.

\section{Future work}
In what follows we provide a sketch of the topics that will be dealt with and will constitute the natural continuation of the work done so far. Our main objective will be, on the basis of the results obtained until now, and by iterating them, to define and investigate the category of towers, whose objects will consist of families (indexed by $\mathbb{N}$) of higher-order many-sorted rewriting systems and of higher-order many-sorted categorial algebras, including higher-order Curry-Howard type results for the latter, together with an additional structure that intertwines such $\omega$-families; and whose morphism from a tower to another will be morphisms between its higher-order many-sorted rewriting systems and higher-order many-sorted categorial algebras compatible with their structures. 

\begin{wrapfigure}{r}{0.45\textwidth}
\vspace{-\baselineskip}
\vspace{-\baselineskip}
\begin{tikzpicture}
[ACliment/.style={-{To [angle'=45, length=5.75pt, width=4pt, round]}},
RHACliment/.style={right hook-{To [angle'=45, length=5.75pt, width=4pt, round]}, font=\scriptsize}, 
scale=1
]

\tikzset{encercla/.style={draw=black, line width=.5pt, inner sep=0pt, rectangle, rounded corners}};

\node[] (A) at (4,0) {$\mathbf{T}_{\Sigma}(X)$};

\node[] (Ap) at (0,0) {$X$};
\draw[ACliment]  (Ap) to node [above] {$\eta^{X}$} (A);

\node[] (B) at (4,-3) {
\begin{tabular}{c}
$\mathbf{T}_{\boldsymbol{\mathcal{E}}^{\boldsymbol{\mathcal{A}}}}(\mathbf{Pth}_{\boldsymbol{\mathcal{A}}})$
\end{tabular}
};
\node[] (Bp) at (0,-3) {$\mathcal{A}$};

\draw[ACliment]  ($(A)+(0,-.2)$) to node [midway, fill=white, rotate=90] {
$\mathrm{ip}^{(\llbracket 1\rrbracket, \llbracket 0\rrbracket)\sharp}$} ($(B)+(0,.2)$);
\draw[ACliment]  ($(B)+(.3,.2)$) to node [below, fill=white, rotate=90] {
$\mathrm{tg}^{(\llbracket 0\rrbracket, \llbracket 1\rrbracket)}$} ($(A)+(.3,-.2)$);
\draw[ACliment]  ($(B)+(-.3,.2)$) to node [above, fill=white, rotate=90] {
$\mathrm{sc}^{(\llbracket 0\rrbracket, \llbracket 1\rrbracket)}$} ($(A)+(-.3,-.2)$);

\draw[ACliment]  (Bp) to node [below, fill=white] {
$\mathrm{ech}^{\mathcal{A}}$} (B);
\draw[ACliment, rotate=30]  ($(Bp)+(.3,.15)$) to node [left] {
$\pi^{\mathcal{A}}_{0}$} ($(A)+(-1,.15)$);
\draw[ACliment, rotate=30]  ($(Bp)+(.3,-.15)$) to node  [right] {
$\pi^{\mathcal{A}}_{1}$} ($(A)+(-1,-.15)$);

\node[] (C) at (4,-6) {
\begin{tabular}{c}
$\mathbf{T}_{\boldsymbol{\mathcal{E}}^{\boldsymbol{\mathcal{A}}^{(2)}}}(\mathbf{Pth}_{\boldsymbol{\mathcal{A}}^{(2)}})$
\end{tabular}
};
\node[] (Cp) at (0,-6) {$\mathcal{A}^{(2)}$};

\draw[ACliment]  ($(B)+(0,-.2)$) to node [midway, fill=white, rotate=90] {
$\mathrm{ip}^{(\llbracket 2\rrbracket, \llbracket 1\rrbracket)\sharp}$} ($(C)+(0,.2)$);
\draw[ACliment]  ($(C)+(.3,.2)$) to node [below, fill=white, rotate=90] {
$\mathrm{tg}^{(\llbracket 1\rrbracket, \llbracket 2\rrbracket)}$} ($(B)+(.3,-.2)$);
\draw[ACliment]  ($(C)+(-.3,.2)$) to node [above, fill=white, rotate=90] {
$\mathrm{sc}^{(\llbracket 1\rrbracket, \llbracket 2\rrbracket)}$} ($(B)+(-.3,-.2)$);

\draw[ACliment]  (Cp) to node [below, fill=white] {
$\mathrm{ech}^{\mathcal{A}^{(2)}}$} (C);
\draw[ACliment, rotate=30]  ($(Cp)+(.3,.15)$) to node [left] {
$\pi^{\mathcal{A}^{(2)}}_{0}$} ($(B)+(-1,.15)$);
\draw[ACliment, rotate=30]  ($(Cp)+(.3,-.15)$) to node  [right] {
$\pi^{\mathcal{A}^{(2)}}_{1}$} ($(B)+(-1,-.15)$);

\node[rotate=270] () at (2,-7.5) {$\cdots$};

\node[] (D) at (4,-9) {
\begin{tabular}{c}
$\mathbf{T}_{\boldsymbol{\mathcal{E}}^{\boldsymbol{\mathcal{A}}^{(n-1)}}}(\mathbf{Pth}_{\boldsymbol{\mathcal{A}}^{(n-1)}})$
\end{tabular}
};

\node[] (E) at (4,-12) {
\begin{tabular}{c}
$\mathbf{T}_{\boldsymbol{\mathcal{E}}^{\boldsymbol{\mathcal{A}}^{(n)}}}(\mathbf{Pth}_{\boldsymbol{\mathcal{A}}^{(n)}})$
\end{tabular}
};
\node[] (Ep) at (0,-12) {$\mathcal{A}^{(n)}$};

\draw[ACliment]  ($(D)+(0,-.2)$) to node [midway, pos=.46, fill=white, rotate=90] {
$\mathrm{ip}^{(\llbracket n\rrbracket, \llbracket n-1\rrbracket)\sharp}$} ($(E)+(0,.3)$);
\draw[ACliment]  ($(E)+(.3,.2)$) to node [below, fill=white, rotate=90] {
$\mathrm{tg}^{(\llbracket n-1\rrbracket, \llbracket n\rrbracket)}$} ($(D)+(.3,-.2)$);
\draw[ACliment]  ($(E)+(-.3,.2)$) to node [above, fill=white, rotate=90] {
$\mathrm{sc}^{(\llbracket n-1\rrbracket, \llbracket n\rrbracket)}$} ($(D)+(-.3,-.2)$);

\draw[ACliment]  (Ep) to node [below, fill=white] {
$\mathrm{ech}^{\mathcal{A}^{(n)}}$} (E);
\draw[ACliment, rotate=30]  ($(Ep)+(.3,.15)$) to node [left] {
$\pi^{\mathcal{A}^{(n)}}_{0}$} ($(D)+(-1,.15)$);
\draw[ACliment, rotate=30]  ($(Ep)+(.3,-.15)$) to node  [right] {
$\pi^{\mathcal{A}^{(n)}}_{1}$} ($(D)+(-1,-.15)$);

\node[rotate=270] () at (2,-14) {$\cdots$};
\end{tikzpicture}
\caption{The tower associated to a rewriting system.}%
\label{FA1}
\vspace{-\baselineskip}
\end{wrapfigure}

Our \textsf{first objective} will be, for a given zeroth order $S$-sorted rewriting system $(S,\Sigma, X)$, to define, by a multiple recursion, families  
$\boldsymbol{\mathcal{A}}^{(\bigcdot)} = (\boldsymbol{\mathcal{A}}^{(n)})_{n\in \mathbb{N}}$, of higher-order many-sorted rewriting systems, where $\boldsymbol{\mathcal{A}}^{(0)}=(S,\Sigma, X)$ and, for $n\in \mathbb{N}$, $\boldsymbol{\mathcal{A}}^{(n+1)}= (\boldsymbol{\mathcal{A}}^{(n)},\mathcal{A}^{(n+1)})$, where $\boldsymbol{\mathcal{A}}^{(n)}$ is an $n$-th order $S$-sorted rewriting system and $\mathcal{A}^{(n+1)}$, the rewriting rules of $(n+1)$-th order, is a choice of a suitable subset of pairs in $\mathrm{T}_{\boldsymbol{\mathcal{E}}^{\boldsymbol{\mathcal{A}}^{(n)}}}(\mathbf{Pth}_{\boldsymbol{\mathcal{A}}^{(n)}})$, the underlying $S$-sorted set of $\mathbf{T}_{\boldsymbol{\mathcal{E}}^{\boldsymbol{\mathcal{A}}^{(n)}}}(\mathbf{Pth}_{\boldsymbol{\mathcal{A}}^{(n)}})$, the free partial $\Sigma^{\boldsymbol{\mathcal{A}}^{(n)}}$-algebra in the category $\mathsf{PAlg}(\boldsymbol{\mathcal{E}}^{\boldsymbol{\mathcal{A}}^{(n)}})$, of partial algebras (see below), generated by the partial 
$\Sigma^{\boldsymbol{\mathcal{A}}^{(n)}}$-algebra $\mathbf{Pth}_{\boldsymbol{\mathcal{A}}^{(n)}}$, of paths of $n$-th order relative to $\boldsymbol{\mathcal{A}}^{(n)}$.  In this regard, let us note that $\Sigma^{\boldsymbol{\mathcal{A}}^{(\bigcdot)}} = (\Sigma^{\boldsymbol{\mathcal{A}}^{(n)}})_{n\in \mathbb{N}}$ is a family of $S$-sorted categorial signatures associated to 
$\boldsymbol{\mathcal{A}}^{(\bigcdot)}$, i.e., $S$-sorted signatures of increasing complexity, where  $\Sigma^{\boldsymbol{\mathcal{A}}^{(0)}}=\Sigma$ and, for $n\in \mathbb{N}$, $\Sigma^{\boldsymbol{\mathcal{A}}^{(n+1)}}$ is the $S$-sorted signature that enlarges $\Sigma^{\boldsymbol{\mathcal{A}}^{(n)}}$ by adding the $n$-th categorial operations of $n$-source, $n$-target and $n$-composition for every sort, and constant operation symbols for every $(n+1)$-th order rewriting rule of every sort in $\mathcal{A}^{(n+1)}$. Note that for simplicity and for $n\geq 1$ we have omitted the explicit dependence on the many-sorted signature $(S,\Sigma^{\boldsymbol{\mathcal{A}}^{(n)}})$ in $\mathbcal{A}^{(n)}$. Moreover, 
$\mathbf{Pth}_{\mathbcal{A}^{(\bigcdot)}} = (\mathbf{Pth}_{\boldsymbol{\mathcal{A}}^{(n)}})_{n\in \mathbb{N}}$ is a family of partial 
$\Sigma^{\boldsymbol{\mathcal{A}}^{(n)}}$-algebra, where $\mathbf{Pth}_{\boldsymbol{\mathcal{A}}^{(0)}}$ is $\mathbf{T}_{\Sigma}(X)$, the free $\Sigma$-algebra on $X$, and, for $n\geq 1$, $\mathbf{Pth}_{\boldsymbol{\mathcal{A}}^{(n)}}$ is the partial $\Sigma^{\boldsymbol{\mathcal{A}}^{(n)}}$-algebra of paths of 
$n$-th order relative to $\boldsymbol{\mathcal{A}}^{(n)}$; and $\mathbf{Pth}_{\mathbcal{A}^{(\bigcdot)}}$ has associated a family $(\coprod \mathrm{Pth}_{\mathbcal{A}^{(n)}},\leq_{\mathbf{Pth}_{\mathbcal{A}^{(n)}}})_{n\in \mathbb{N}}$ of Artinian preordered sets which will allow us to make proofs by  Artinian induction and definitions by Artinian recursion at points that are key to the development of this work. 
Furthermore, $\mathsf{PAlg}(\boldsymbol{\mathcal{E}}^{(\boldsymbol{\mathcal{A}}^{(\bigcdot)})}) = (\mathsf{PAlg}(\boldsymbol{\mathcal{E}}^{(\boldsymbol{\mathcal{A}}^{(n)})}))_{n\in \mathbb{N}}$ is a family of categories of partial $\Sigma^{\boldsymbol{\mathcal{A}}^{(\bigcdot)}}$-algebras, where 
$\mathsf{PAlg}(\boldsymbol{\mathcal{E}}^{(\boldsymbol{\mathcal{A}}^{(0)})})$ is the category of many-sorted total $\Sigma^{\boldsymbol{\mathcal{A}}^{(0)}}$-algebras, i.e., the category $\mathsf{Alg}(\Sigma)$, and, for $n\geq 1$, $\mathsf{PAlg}(\boldsymbol{\mathcal{E}}^{(\boldsymbol{\mathcal{A}}^{(n)})})$ is the category canonically associated to the $\mathrm{QE}$-variety $\mathcal{V}(\boldsymbol{\mathcal{E}}^{(\boldsymbol{\mathcal{A}}^{(n)})})$, of partial $\Sigma^{\boldsymbol{\mathcal{A}}^{(n)}}$-algebras, determined by the specification $\boldsymbol{\mathcal{E}}^{(\boldsymbol{\mathcal{A}}^{(n)})}$ associated to $\boldsymbol{\mathcal{A}}^{(n)}$, whose underlying system of defining equations $\mathcal{E}^{(\boldsymbol{\mathcal{A}}^{(n)})}$ is formed by $\mathrm{QE}$-equations. See~\cite{an83, bs67, pb82}, for the notions of $\mathrm{QE}$-variety and $\mathrm{QE}$-equations.
Structures in $\mathsf{PAlg}(\boldsymbol{\mathcal{E}}^{(\boldsymbol{\mathcal{A}}^{(n)})})$, i.e., models of 
$\mathcal{E}^{(\boldsymbol{\mathcal{A}}^{(n)})}$, are instances of what we will call in this work $S$-sorted $n$-categorial $\Sigma$-algebras. Finally, $\mathbf{T}_{\boldsymbol{\mathcal{E}}^{\boldsymbol{\mathcal{A}}^{(\bigcdot)}}}(\mathbf{Pth}_{\boldsymbol{\mathcal{A}}^{(\bigcdot)}})=(\mathbf{T}_{\boldsymbol{\mathcal{E}}^{\boldsymbol{\mathcal{A}}^{(n)}}}(\mathbf{Pth}_{\boldsymbol{\mathcal{A}}^{(n)}}))_{n\in \mathbb{N}}$, as we have said above, is the family of the free partial $\Sigma^{\boldsymbol{\mathcal{A}}^{(\bigcdot)}}$-algebras generated by $\mathbf{Pth}_{\boldsymbol{\mathcal{A}}^{(\bigcdot)}}$.

The  involved structures are also intertwined through different mappings. For each $n\in \mathbb{N}$, since $\mathcal{A}^{(n+1)}$, the rewriting rules of $(n+1)$-th order, is a choice of a suitable subset of pairs in $\mathrm{T}_{\boldsymbol{\mathcal{E}}^{\boldsymbol{\mathcal{A}}^{(n)}}}(\mathbf{Pth}_{\boldsymbol{\mathcal{A}}^{(n)}})$, we will have, for $i\in 2$, the canonical projections, denoted by $\pi^{\mathcal{A}^{(n+1)}}_{i}$, from $\mathcal{A}^{(n+1)}$ to $\mathrm{T}_{\boldsymbol{\mathcal{E}}^{\boldsymbol{\mathcal{A}}^{(n)}}}(\mathbf{Pth}_{\boldsymbol{\mathcal{A}}^{(n)}})$. Moreover, for every $n\in \mathbb{N}$ and every sort $s\in S$, there exists a canonical interpretation of a rewriting rule of sort $s$ in $\mathcal{A}^{(n+1)}$ as a path of length $1$, what in this work has been called an $(n+1)$-th order echelon. This interpretation is the one that serves, for example, to describe the underlying complexity of a path, which is, the foundation for defining the partial Artinian preorder $(\coprod \mathrm{Pth}_{\mathbcal{A}^{(n)}},\leq_{\mathbf{Pth}_{\mathbcal{A}^{(n)}}})_{n\in \mathbb{N}}$. Therefore, the family of echelon interpretations $\mathrm{ech}^{\mathcal{A}^{(\bigcdot)}} = (\mathrm{ech}^{\mathcal{A}^{(n)}})_{n\in \mathbb{N}}$ is a family of $S$-sorted mappings from
$\mathcal{A}^{(\bigcdot)}$ to $\mathrm{T}_{\mathbcal{E}^{\boldsymbol{\mathcal{A}}^{(\bigcdot)}}}(\mathbf{Pth}_{\mathbcal{A}^{(\bigcdot)}})$, where  $\mathrm{ech}^{\mathcal{A}^{(0)}}=\eta^{X}$, the canonical insertion of generators, from $X$ to $\mathrm{T}_{\Sigma}(X)$,  and, for $n\in \mathbb{N}$, $\mathrm{ech}^{\mathcal{A}^{(n+1)}}$ is the echelon interpretation of $(n+1)$-th order rewriting rules in $\mathcal{A}^{(n+1)}$ as $(n+1)$-th order echelons in $\mathrm{T}_{\mathbcal{E}^{\boldsymbol{\mathcal{A}}^{(n+1)}}}(\mathbf{Pth}_{\mathbcal{A}^{(n+1)}})$. Additionally, for each $n\in \mathbb{N}$, the layers $\mathrm{T}_{\mathbcal{E}^{\boldsymbol{\mathcal{A}}^{(n+1)}}}(\mathbf{Pth}_{\mathbcal{A}^{(n+1)}})$ and $\mathrm{T}_{\mathbcal{E}^{\boldsymbol{\mathcal{A}}^{(n)}}}(\mathbf{Pth}_{\mathbcal{A}^{(n)}})$ are intertwined with two downstream applications, $\mathrm{sc}^{(\llbracket n\rrbracket, \llbracket n+1\rrbracket)}$ and $\mathrm{tg}^{(\llbracket n\rrbracket, \llbracket n+1\rrbracket)}$, from $\mathrm{T}_{\mathbcal{E}^{\boldsymbol{\mathcal{A}}^{(n+1)}}}(\mathbf{Pth}_{\mathbcal{A}^{(n+1)}})$ to $\mathrm{T}_{\mathbcal{E}^{\boldsymbol{\mathcal{A}}^{(n)}}}(\mathbf{Pth}_{\mathbcal{A}^{(n)}})$, sending a $(n+1)$-th order path to its $n$-th order source and target, respectively, and a upstream application $\mathrm{ip}^{(\llbracket n+1\rrbracket, \llbracket n\rrbracket)\sharp}$, from $\mathrm{T}_{\mathbcal{E}^{\boldsymbol{\mathcal{A}}^{(n)}}}(\mathbf{Pth}_{\mathbcal{A}^{(n)}})$ to $\mathrm{T}_{\mathbcal{E}^{\boldsymbol{\mathcal{A}}^{(n+1)}}}(\mathbf{Pth}_{\mathbcal{A}^{(n+1)}})$, sending an $n$-th order path to the $(n+1)$-th order identity path on it. These latter mappings are obtained by means of the universal property of $\mathrm{T}_{\mathbcal{E}^{\boldsymbol{\mathcal{A}}^{(n)}}}(\mathbf{Pth}_{\mathbcal{A}^{(n)}})$. We will show that these mappings are, precisely, $\Sigma^{\boldsymbol{\mathcal{A}}^{(n)}}$-homomorphisms between the $\Sigma^{\boldsymbol{\mathcal{A}}^{(n)}}$-algebras $\mathbf{T}_{\mathbcal{E}^{\boldsymbol{\mathcal{A}}^{(n)}}}(\mathbf{Pth}_{\mathbcal{A}^{(n)}})$ and $\mathbf{T}^{(n,n+1)}_{\mathbcal{E}^{\boldsymbol{\mathcal{A}}^{(n+1)}}}(\mathbf{Pth}_{\mathbcal{A}^{(n+1)}})$, the $\Sigma^{\boldsymbol{\mathcal{A}}^{(n)}}$-reduction of the  $\Sigma^{\boldsymbol{\mathcal{A}}^{(n+1)}}$-algebra $\mathbf{T}_{\mathbcal{E}^{\boldsymbol{\mathcal{A}}^{(n+1)}}}(\mathbf{Pth}_{\mathbcal{A}^{(n+1)}})$.

%

The $7$-tuple
\[
\begin{adjustbox}{width=\textwidth}$
\left(
\mathbcal{A}^{(\bigcdot)},
\mathbf{T}_{\mathbcal{E}^{\boldsymbol{\mathcal{A}}^{(\bigcdot)}}}(\mathbf{Pth}_{\mathbcal{A}^{(\bigcdot)}}),
(\pi^{\mathcal{A}^{(\bigcdot^{+})}}_{i})_{i\in 2},
\mathrm{ech}^{\mathcal{A}^{(\bigcdot)}},
\mathrm{sc}^{(\llbracket (\bigcdot)\rrbracket, \llbracket (\bigcdot^{+})\rrbracket)},
\mathrm{tg}^{(\llbracket (\bigcdot)\rrbracket, \llbracket (\bigcdot^{+})\rrbracket)},
\mathrm{ip}^{(\llbracket (\bigcdot^{+})\rrbracket, \llbracket (\bigcdot)\rrbracket)\sharp}
\right)$
\end{adjustbox}
\]
denoted simply by $(\mathbcal{A}^{(\bigcdot)},\mathbf{T}_{\mathbcal{E}^{\boldsymbol{\mathcal{A}}^{(\bigcdot)}}}(\mathbf{Pth}_{\mathbcal{A}^{(\bigcdot)}}))$ or by $\mathbb{A}$, for short, will be called a tower (metaphorically speaking, $\boldsymbol{\mathcal{A}}^{(\bigcdot)}$ is the seed the tower $\mathbb{A}$ grows up from). An $\mathbcal{A}^{(\bigcdot)}$-tower will be graphically represented, in full, as shown in the diagram in Figure~\ref{FA1}.

Regarding the family $\mathbf{T}_{\boldsymbol{\mathcal{E}}^{\boldsymbol{\mathcal{A}}^{(\bigcdot)}}}(\mathbf{Pth}_{\boldsymbol{\mathcal{A}}^{(\bigcdot)}})$, of free partial $\Sigma^{\boldsymbol{\mathcal{A}}^{(\bigcdot)}}$-algebras, we provide two explicit isomorphic representations of it. The first explicit representation, namely $\llbracket \mathbf{Pth}_{\boldsymbol{\mathcal{A}}^{(\bigcdot)}}\rrbracket=(\llbracket \mathbf{Pth}_{\boldsymbol{\mathcal{A}}^{(n)}}\rrbracket)_{n\in\mathbb{N}}$, is a quotient of the family $\mathbf{Pth}_{\boldsymbol{\mathcal{A}}^{(\bigcdot)}}$ of paths, i.e., for every $n\in\mathbb{N}$, $\llbracket \mathbf{Pth}_{\boldsymbol{\mathcal{A}}^{(n)}}\rrbracket$ is a quotient of $\mathbf{Pth}_{\boldsymbol{\mathcal{A}}^{(n)}}$.  Furthermore, taking into account that, for every $n\in\mathbb{N}$, $\llbracket \mathbf{Pth}_{\boldsymbol{\mathcal{A}}^{(n)}}\rrbracket$ is an $S$-sorted $n$-categorial $\Sigma$-algebra in the $\mathrm{QE}$-variety $\mathcal{V}(\boldsymbol{\mathcal{E}}^{(\boldsymbol{\mathcal{A}}^{(n)})})$ that is $\mathbf{Pth}_{\boldsymbol{\mathcal{A}}^{(n)}}$-generated, by the universal property of $\mathbf{T}_{\boldsymbol{\mathcal{E}}^{\boldsymbol{\mathcal{A}}^{(n)}}}(\mathbf{Pth}_{\boldsymbol{\mathcal{A}}^{(n)}})$, the quotient mapping $\mathrm{pr}^{\llbracket n\rrbracket}$, from $\mathbf{Pth}_{\boldsymbol{\mathcal{A}}^{(n)}}$ to $\llbracket \mathbf{Pth}_{\boldsymbol{\mathcal{A}}^{(n)}}\rrbracket$ can be extended to a unique $\Sigma^{\boldsymbol{\mathcal{A}}^{(n)}}$-homomorphism, which we write by $\mathrm{pr}^{\llbracket n\rrbracket\mathsf{p}}$, from $\mathbf{T}_{\boldsymbol{\mathcal{E}}^{\boldsymbol{\mathcal{A}}^{(n)}}}(\mathbf{Pth}_{\boldsymbol{\mathcal{A}}^{(n)}})$ to $\llbracket \mathbf{Pth}_{\boldsymbol{\mathcal{A}}^{(n)}}\rrbracket$, and which will be a $\Sigma^{\boldsymbol{\mathcal{A}}^{(\bigcdot)}}$-isomorphism.

The second explicit representation of $\mathbf{T}_{\boldsymbol{\mathcal{E}}^{\boldsymbol{\mathcal{A}}^{(n)}}}(\mathbf{Pth}_{\boldsymbol{\mathcal{A}}^{(n)}})$, namely $\llbracket \mathbf{PT}_{\boldsymbol{\mathcal{A}}^{(\bigcdot)}}\rrbracket=(\llbracket \mathbf{PT}_{\boldsymbol{\mathcal{A}}^{(n)}}\rrbracket)_{n\in\mathbb{N}}$, the family of path-term classes, is a subquotient of the family $\mathbf{T}_{\Sigma^{\boldsymbol{\mathcal{A}}^{(\bigcdot)}}}(X)=(\mathbf{T}_{\Sigma^{\boldsymbol{\mathcal{A}}^{(n)}}}(X))_{n\in\mathbb{N}}$, of free $\Sigma^{\boldsymbol{\mathcal{A}}^{(\bigcdot)}}$-algebras, i.e., for every $n\in\mathbb{N}$, $\llbracket \mathbf{PT}_{\boldsymbol{\mathcal{A}}^{(n)}}\rrbracket$ is a subquotient of the free $\Sigma^{\boldsymbol{\mathcal{A}}^{(n)}}$-algebra $\mathbf{T}_{\Sigma^{\boldsymbol{\mathcal{A}}^{(n)}}}(X)$. To explicitly construct the aforementioned isomorphism between the families $\llbracket \mathbf{Pth}_{\boldsymbol{\mathcal{A}}^{(\bigcdot)}}\rrbracket$ and $\llbracket \mathbf{PT}_{\boldsymbol{\mathcal{A}}^{(\bigcdot)}}\rrbracket$ we rely on the fact that $\mathbf{Pth}_{\boldsymbol{\mathcal{A}}^{(\bigcdot)}}$ admits an Artinian preorder compatible with the quotient structure on $\llbracket \mathbf{Pth}_{\boldsymbol{\mathcal{A}}^{(\bigcdot)}}\rrbracket$. This will allow us to carry out proofs by Artinian induction and definitions by Artinian recursion on the families of paths and quotients of paths.
Applying what has just been said, we explicitly construct a family $\mathrm{CH}^{\llbracket\bigcdot\rrbracket}=(\mathrm{CH}^{\llbracket n\rrbracket})_{n\in\mathbb{N}}$ of monotone $\Sigma^{\boldsymbol{\mathcal{A}}^{(\bigcdot)}}$-isomorphisms from $\llbracket \mathbf{Pth}_{\boldsymbol{\mathcal{A}}^{(\bigcdot)}}\rrbracket$ to $\llbracket \mathbf{PT}_{\boldsymbol{\mathcal{A}}^{(\bigcdot)}}\rrbracket$, which we call the Curry-Howard isomorphisms, i.e., for every $n\in\mathbb{N}$, we provide $\mathrm{CH}^{\llbracket n\rrbracket}$, a Curry-Howard monotone $\Sigma^{\boldsymbol{\mathcal{A}}^{(n)}}$-isomorphism from  $\llbracket \mathbf{Pth}_{\boldsymbol{\mathcal{A}}^{(n)}}\rrbracket$ to $\llbracket \mathbf{PT}_{\boldsymbol{\mathcal{A}}^{(n)}}\rrbracket$.  

These isomorphisms, together with their inverses, are instances of what we will call in this work $S$-sorted $n$-categorial $\Sigma$-homomorphisms between $S$-sorted $n$-categorial $\Sigma$-algebras, see Figure~\ref{FA2}.

\begin{figure}
\begin{tikzpicture}
[ACliment/.style={-{To [angle'=45, length=5.75pt, width=4pt, round]}},
RHACliment/.style={right hook-{To [angle'=45, length=5.75pt, width=4pt, round]}, font=\scriptsize}, 
scale=.9
]

\tikzset{encercla/.style={draw=black, line width=.5pt, inner sep=0pt, rectangle, rounded corners}};

\node[] (Q1) at (-1.6,-1) {};
\node[] (Q2) at (9.5,1) {};

\node[] (Anp) at (4.5,0) {$\llbracket \mathbf{Pth}_{\boldsymbol{\mathcal{A}}^{(n)}}\rrbracket$};
\node[] (Anpt) at (8.5,0) {$\llbracket \mathbf{PT}_{\boldsymbol{\mathcal{A}}^{(n)}}\rrbracket$};
\node[] (Anptt) at (0,0) {$\mathbf{T}_{\boldsymbol{\mathcal{E}}^{\boldsymbol{\mathcal{A}}^{(n)}}}(\mathbf{Pth}_{\boldsymbol{\mathcal{A}}^{(n)}})$};


\node[] (iso1) at (2.5,0) {$\cong$};
\draw[ACliment]  ($(iso1)+(-1,.3)$) to node [ above] {
$\textstyle \mathrm{pr}^{\llbracket n\rrbracket\mathsf{p}}$
} ($(iso1)+(1,.3)$);
\draw[ACliment]  ($(iso1)+(1,-.3)$) to node [below] {
} ($(iso1)+(-1,-.3)$);

\node[] (iso2) at (6.5,0) {$\cong$};
\draw[ACliment]  ($(iso2)+(-1,.3)$) to node [ above] {
$\textstyle \mathrm{CH}^{\llbracket n\rrbracket}$} ($(iso2)+(1,.3)$);
\draw[ACliment]  ($(iso2)+(1,-.3)$) to node [below] {
} ($(iso2)+(-1,-.3)$);

\end{tikzpicture}
\caption{Isomorphisms of $S$-sorted $n$-categorial $\Sigma$-algebras.}
\label{FA2}
\end{figure}

Towers, like any other kind of species of structures, are nothing more than the static component of a whole whose dynamic component, following the category-theoretic dictum, is given by morphisms between towers.
 
Our \textsf{second objective} will therefore be, given towers $\mathbb{A} = (\mathbcal{A}^{(\bigcdot)},\mathbf{T}_{\mathbcal{E}^{\boldsymbol{\mathcal{A}}^{(\bigcdot)}}}(\mathbf{Pth}_{\mathbcal{A}^{(\bigcdot)}}))$ and $\mathbb{B} = (\mathbcal{B}^{(\bigcdot)},\mathbf{T}_{\mathbcal{E}^{\boldsymbol{\mathcal{B}}^{(\bigcdot)}}}(\mathbf{Pth}_{\mathbcal{B}^{(\bigcdot)}}))$ to define a notion of morphism from $\mathbb{A}$ to $\mathbb{B}$ in such a way that (after defining the identity morphisms, the composition and verifying the category axioms) we get a category. Specifically, a morphism from $\mathbb{A}$ to $\mathbb{B}$ will be an
ordered triple $(\mathbb{A},\mathbf{f}^{\llbracket\bigcdot\rrbracket},\mathbb{B})$, denoted by 
$\mathbf{f}^{\llbracket\bigcdot\rrbracket}\colon \mathbb{A}\mor \mathbb{B}$ or $\mathbf{f}^{\llbracket \bigcdot \rrbracket}$ 
for short, in which $\mathbf{f}^{\llbracket\bigcdot\rrbracket}$ is a family $(\mathbf{f}^{\llbracket n \rrbracket@})_{n\in \mathbb{N}}$ such that $\mathbf{f}^{\llbracket 0 \rrbracket @}$ will be defined explicitly and, for every $n\in \mathbb{N}$, $\mathbf{f}^{\llbracket n+1 \rrbracket@} = (\llbracket \mathbf{f}^{(n+1)} \rrbracket,f^{\llbracket n+1 \rrbracket @})$ where 
$\llbracket \mathbf{f}^{(n+1)} \rrbracket$ is a $\cong^{(n+1)}$-class of a morphism from $\mathbcal{A}^{(n+1)}$ to $\mathbcal{B}^{(n+1)}$ and $f^{\llbracket n+1 \rrbracket@}$ a homomorphism from $\mathbf{T}_{\mathbcal{E}^{\boldsymbol{\mathcal{A}}^{(n+1)}}}(\mathbf{Pth}_{\mathbcal{A}^{(n+1)}})$ to $\mathbf{T}_{\mathbcal{E}^{\boldsymbol{\mathcal{B}}^{(n+1)}}}(\mathbf{Pth}_{\mathbcal{B}^{(n+1)}})_{\varphi}$, compatible with the structural morphisms of the involved towers, and defined by simultaneous recursion as follows. 


For $n = 0$ we have that $\mathbcal{A}^{(0)} = (S,\Sigma,X)$ and $\mathbcal{B}^{(0)} = (T,\Lambda,Y)$ are zeroth order many-sorted rewriting systems, and that a morphism from the tower $\mathbb{A}^{(0)}$ to the tower $\mathbb{B}^{(0)}$ will be an ordered triple 
$(\mathbb{A}^{(0)},\mathbf{f}^{(0)\sharp},\mathbb{B}^{(0)})$, denoted by 
$\mathbf{f}^{(0)\sharp}\colon \mathbb{A}^{(0)}\mor \mathbb{B}^{(0)}$ or $\mathbf{f}^{(0)}$ 
for short, in which $\mathbf{f}^{(0)\sharp}$ is an ordered pair $(\mathbf{f}^{(0)},f^{(0)\sharp})$ where $\mathbf{f}^{(0)} = (\varphi, c, f^{(0)})$ is a morphism from $(S,\Sigma,X)$ to $(T,\Lambda,Y)$ (i.e., a morphism of $\mathsf{Rws}_{\mathfrak{d}}^{(0)}$) and, since $\mathbf{T}_{\mathbcal{E}^{\boldsymbol{\mathcal{A}}^{(0)}}}(\mathbf{Pth}_{\mathbcal{A}^{(0)}})$ is $\mathbf{T}_{\Sigma}(X)$ and $\mathbf{T}_{\mathbcal{E}^{\boldsymbol{\mathcal{B}}^{(0)}}}(\mathbf{Pth}_{\mathbcal{B}^{(0)}})$ is $\mathbf{T}_{\Lambda}(Y)$, $f^{(0)\sharp}$ is the unique $\Sigma$-homomorphism from $\mathbf{T}_{\Sigma}(X)$ to $\mathbf{c}^{\ast}_{\mathfrak{d}}(\mathbf{T}_{\Lambda}(Y))$ such that $f^{(0)} = f^{(0)\sharp}\circ \eta^{X}$. 


Let us suppose defined $\mathbf{f}^{\llbracket n \rrbracket@} = (\llbracket\mathbf{f}^{(n)}\rrbracket,f^{\llbracket n \rrbracket@})$ from 
$\mathbb{A}^{(n)}$ to $\mathbb{B}^{(n)}$, for $n\geq 1$. Then we have that $\mathbb{A}^{(n+1)}$ and $\mathbb{B}^{(n+1)}$ are $(n+1)$-th order towers associated with the $(n+1)$-th order many-sorted rewriting systems $\mathbcal{A}^{(n+1)}$ and $\mathbcal{B}^{(n+1)}$, respectively, and that a morphism from $\mathbb{A}^{(n+1)}$ to $\mathbb{B}^{(n+1)}$ will be an ordered triple 
$(\mathbb{A}^{(n+1)},\mathbf{f}^{\llbracket n+1 \rrbracket @},\mathbb{B}^{(n+1)}),$ 
denoted by $\mathbf{f}^{\llbracket n+1 \rrbracket@}\colon \mathbcal{A}^{(n+1)}\mor \mathbcal{B}^{(n+1)}$ or $\mathbf{f}^{\llbracket n+1 \rrbracket@}$ for short, in which $\mathbf{f}^{\llbracket n+1 \rrbracket@}$ is an ordered pair  $(\llbracket \mathbf{f}^{(n+1)} \rrbracket,f^{\llbracket n+1 \rrbracket@})$ where
\begin{enumerate}
\item $\llbracket \mathbf{f}^{(n+1)} \rrbracket$ is a $\cong^{(n+1)}$-class of a morphism from $\mathbcal{A}^{(n+1)}$ to $\mathbcal{B}^{(n+1)}$, where $\mathbf{f}^{(n+1)} = (\mathbf{f}^{(n)}, f^{(n+1)})$ and $f^{(n+1)}$ is an $S$-sorted mapping from $\mathcal{A}^{(n+1)}$ to $\mathrm{Pth}_{\mathbcal{B}^{(n+1)},\varphi}$ such that, for every $s\in S$ and every $\mathfrak{p}^{(n+1)}\in \mathcal{A}^{(n+1)}_{s}$, 
we have that 
$$
f^{(n+1)}_{s}(\mathfrak{p}^{(n+1)})\in \mathrm{Pth}_{\mathbcal{B}^{(n+1)},\varphi(s)}
(
f^{\llbracket n \rrbracket@}_{s}(
\pi_{0,s}^{\mathcal{A}^{(n+1)}}(\mathfrak{p}^{(n+1)})
),
f^{\llbracket n \rrbracket@}_{s}(
\pi_{1,s}^{\mathcal{A}^{(n+1)}}(\mathfrak{p}^{(n+1)})
)
),
$$
and
\item $f^{\llbracket n+1 \rrbracket@}$ is the homomorphism univocally determined by $\llbracket \mathbf{f}^{(n+1)} \rrbracket$.
\end{enumerate}
We let $\mathsf{Rws}_{\mathfrak{d}}^{\llbracket n+1 \rrbracket}$ stand for the category of $(n+1)$-th order many-sorted rewriting systems and $\cong^{(n+1)}$-classes of morphisms between them. 

Let us explain exactly how to obtain $f^{\llbracket n+1 \rrbracket@}$. From $f^{(n+1)}$ we obtain, by recursion on $(\coprod \mathrm{Pth}_{\mathbcal{A}^{(n+1)}},\leq_{\mathbf{Pth}_{\mathbcal{A}}^{(n+1)}})$, a $\Sigma$-homomorphism $f^{(n+1)\flat}$ from $\mathbf{Pth}_{\mathcal{A}^{(n+1)}}^{(0,n+1)}$ to $\mathbf{Pth}_{\mathcal{B}^{(n+1)},\varphi}^{(0,n+1)}$ such that $f^{(n+1)\flat}$ coincides with $f^{\llbracket k \rrbracket@}$ when restricted to the $(n+1,\llbracket k \rrbracket)$-identity paths, for $k\in n+1$, and with $f^{(n+1)}$ when restricted to the $(n+1)$-th order echelons, i.e., to the one-step $n+1$-th order paths canonically associated to the rewrite rules on $\mathcal{A}^{(n+1)}$. In turn, from $f^{(n+1)\flat}$ and since it, as will be proved in this work, is compatible with suitable congruences on $\mathbf{Pth}_{\mathbcal{A}^{(n+1)}}$ and $\mathbf{Pth}_{\mathbcal{B}^{(n+1)}}^{\mathbf{f}^{(n+1)}}$, 
we obtain $f^{\llbracket n+1 \rrbracket@}$, the unique $\Sigma^{\boldsymbol{\mathcal{A}}^{(n+1)}}$-homomorphism  from $\llbracket \mathbf{Pth}_{\boldsymbol{\mathcal{A}}^{(n+1)}}\rrbracket$ to $\llbracket \mathbf{Pth}_{\boldsymbol{\mathcal{B}}^{(n+1)}}^{\mathbf{f}^{(n+1)}}\rrbracket$ induced by $f^{(n+1)\flat}$ on passing to the quotients with respect to the just mentioned congruences. Then, taking into account the isomorphisms stated above, we obtain a unique $\Sigma^{\boldsymbol{\mathcal{A}}^{(n+1)}}$-homomorphism, also denoted by $f^{\llbracket n+1 \rrbracket@}$, from $\mathbf{T}_{\mathbcal{E}^{\boldsymbol{\mathcal{A}}^{(n+1)}}}(\mathbf{Pth}_{\mathbcal{A}^{(n+1)}})$ to 
$\mathbf{T}_{\mathbcal{E}^{\boldsymbol{\mathcal{B}}^{(n+1)}}}^{\mathbf{f}^{(n+1)}}(\mathbf{Pth}_{\mathbcal{B}^{(n+1)}})$ making commutative a suitable diagram. Let us note that, for every $s\in S$ and every $\mathfrak{p}^{(n+1)}\in \mathcal{A}^{(n+1)}_{s}$, the structural operation in $\mathbf{T}_{\mathbcal{E}^{\boldsymbol{\mathcal{B}}^{(n+1)}}}(\mathbf{Pth}_{\mathbcal{B}^{(n+1)}})_{\varphi}$ associated to $\mathfrak{p}^{(n+1)}$, which is an operation of biarity $(\lambda,s)$, i.e., an element of 
$\mathrm{T}_{\mathbcal{E}^{\boldsymbol{\mathcal{B}}^{(n+1)}}}(\mathbf{Pth}_{\mathbcal{B}^{(n+1)}})_{\varphi(s)}$, is $f^{(n+1)}_{s}(\mathfrak{p}^{(n+1)})$ (more precisely, it is the image of $f^{(n+1)}_{s}(\mathfrak{p}^{(n+1)})$ under the canonical embedding of $\mathbf{Pth}_{\mathbcal{B}^{(n+1)}}$ into $\mathbf{T}_{\mathbcal{E}^{\boldsymbol{\mathcal{B}}^{(n+1)}}}(\mathbf{Pth}_{\mathbcal{B}^{(n+1)}})$). 

\begin{figure}
\begin{center}
\begin{tikzpicture}
[ACliment/.style={-{To [angle'=45, length=5.75pt, width=4pt, round]}},
RHACliment/.style={right hook-{To [angle'=45, length=5.75pt, width=4pt, round]}, font=\scriptsize}, 
scale=.95
]

\tikzset{encercla/.style={draw=black, line width=.5pt, inner sep=0pt, rectangle, rounded corners}};


\node[color=white] (A0aux) at (0,0){$\mathbf{Pth}_{\boldsymbol{\mathcal{A}}^{(n)}}/{\mathrm{Ker}(\mathrm{CH}^{(n)})}$};
\node[color=blue!95!black] (A0) at (0,0) {$\mathbf{T}_{\Sigma}(X)$};

\node [encercla, draw=blue!95!black,  fit=(A0aux)] {} ;

\node[] (A0aux) at (0,-1.5) {};
\draw[ACliment, color=blue!95!black]  ($(A0aux)+(0,1)$) to node [sloped, rotate=180, midway, fill=white] {
$\textstyle \mathrm{ip}^{([1],0)\sharp}$
} ($(A0aux)+(0,-1)$);
\draw[ACliment, color=blue!95!black]  ($(A0aux)+(.3,-1)$) to node [sloped, below, fill=white] {
$\textstyle \mathrm{tg}^{(0,[1])}$
} ($(A0aux)+(.3,1)$);
\draw[ACliment, color=blue!95!black]  ($(A0aux)+(-.3,-1)$) to node [sloped, above, fill=white] {
$\textstyle \mathrm{sc}^{(0,[1])}$
} ($(A0aux)+(-.3,1)$);

\node[color=white] (A1paux) at (0,-3) {$\mathbf{Pth}_{\boldsymbol{\mathcal{A}}^{(n)}}/{\mathrm{Ker}(\mathrm{CH}^{(n)})}$};
\node[color=white] (A1ptaux) at (0,-5) {$\mathbf{Pth}_{\boldsymbol{\mathcal{A}}^{(n)}}/{\mathrm{Ker}(\mathrm{CH}^{(n)})}$};

\node [encercla, draw=blue!80!black,  fit=(A1paux)(A1ptaux)] {};

\node[color=blue!80!black] (A1p) at (0,-3) {$[\mathbf{Pth}_{\boldsymbol{\mathcal{A}}}]$};
\node[color=blue!80!black] (A1pt) at (0,-5) {$[\mathbf{PT}_{\boldsymbol{\mathcal{A}}}]$};

\node[rotate=90, color=blue!80!black] (A1iso) at (0,-4) {$\cong$};
\draw[ACliment, color=blue!80!black]  ($(A1iso)+(-.3,.7)$) to node [sloped, rotate=180, above] {
$\textstyle \mathrm{CH}^{[1]}$
} ($(A1iso)+(-.3,-.7)$);
\draw[ACliment, color=blue!80!black]  ($(A1iso)+(.3,-.7)$) to node [sloped, below] {
$\textstyle \mathrm{ip}^{([1],X)@}$
} ($(A1iso)+(.3,.7)$);

\node[] (A1aux) at (0,-6.5) {};
\draw[ACliment, color=blue!80!black]  ($(A1aux)+(0,1)$) to node [sloped, rotate=180, midway, fill=white] {
$\textstyle \mathrm{ip}^{(\llbracket 2\rrbracket,[1])\sharp}$
} ($(A1aux)+(0,-1)$);
\draw[ACliment, color=blue!80!black]  ($(A1aux)+(.3,-1)$) to node [sloped, below, fill=white] {
$\textstyle \mathrm{tg}^{([1],\llbracket 2\rrbracket)}$
} ($(A1aux)+(.3,1)$);
\draw[ACliment, color=blue!80!black]  ($(A1aux)+(-.3,-1)$) to node [sloped, above, fill=white] {
$\textstyle \mathrm{sc}^{([1],\llbracket 2\rrbracket)}$
} ($(A1aux)+(-.3,1)$);

\node[color=white] (A2paux) at (0,-8) {$\mathbf{Pth}_{\boldsymbol{\mathcal{A}}^{(n)}}/{\mathrm{Ker}(\mathrm{CH}^{(n)})}$};
\node[color=white] (A2ptaux) at (0,-10) {$\mathbf{Pth}_{\boldsymbol{\mathcal{A}}^{(n)}}/{\mathrm{Ker}(\mathrm{CH}^{(n)})}$};
\node [encercla, draw=blue!65!black,  fit=(A2paux)(A2ptaux)] {};

\node[color=blue!65!black] (A2p) at (0,-8) {$\llbracket \mathbf{Pth}_{\boldsymbol{\mathcal{A}}^{(2)}}\rrbracket$};
\node[color=blue!65!black] (A2pt) at (0,-10) {$\llbracket \mathbf{PT}_{\boldsymbol{\mathcal{A}}^{(2)}}\rrbracket$};

\node[rotate=90, color=blue!65!black] (A2iso) at (0,-9) {$\cong$};
\draw[ACliment, color=blue!65!black]  ($(A2iso)+(-.3,.7)$) to node [sloped, rotate=180, above] {
$\textstyle \mathrm{CH}^{\llbracket 2 \rrbracket}$
} ($(A2iso)+(-.3,-.7)$);
\draw[ACliment, color=blue!65!black]  ($(A2iso)+(.3,-.7)$) to node [sloped, below] {
$\textstyle \mathrm{ip}^{(\llbracket 2 \rrbracket, X)@}$
} ($(A2iso)+(.3,.7)$);

\node[] (A2aux) at (0,-11.5) {};
\draw[ACliment, color=blue!65!black]  ($(A2aux)+(0,1)$) to node [sloped, rotate=180, midway, fill=white] {
$\textstyle \mathrm{ip}^{(\llbracket 3\rrbracket,\llbracket 2\rrbracket)\sharp}$
} ($(A2aux)+(0,-1)$);
\draw[ACliment, color=blue!65!black]  ($(A2aux)+(.3,-1)$) to node [sloped, below, fill=white] {
$\textstyle \mathrm{tg}^{(\llbracket 2\rrbracket,\llbracket 3\llbracket)}$
} ($(A2aux)+(.3,1)$);
\draw[ACliment, color=blue!65!black]  ($(A2aux)+(-.3,-1)$) to node [sloped, above, fill=white] {
$\textstyle \mathrm{sc}^{(\llbracket 2\rrbracket,\llbracket 3\rrbracket)}$
} ($(A2aux)+(-.3,1)$);

\node[] (ca) at (0,-12.85) {$\vdots$};

\node[] (A3aux) at (0,-14.5) {};
\draw[ACliment, color=blue!50!black]  ($(A3aux)+(0,1)$) to node [sloped, rotate=180, midway, fill=white] {
$\textstyle \mathrm{ip}^{(\llbracket n\rrbracket,\llbracket n-1\rrbracket)\sharp}$
} ($(A3aux)+(0,-1)$);
\draw[ACliment, color=blue!50!black]  ($(A3aux)+(.3,-1)$) to node [sloped, below, fill=white] {
$\textstyle \mathrm{tg}^{(\llbracket n-1\rrbracket,\llbracket n\rrbracket)}$
} ($(A3aux)+(.3,1)$);
\draw[ACliment, color=blue!50!black]  ($(A3aux)+(-.3,-1)$) to node [sloped, above, fill=white] {
$\textstyle \mathrm{sc}^{(\llbracket n-1\rrbracket, \llbracket n\rrbracket)}$
} ($(A3aux)+(-.3,1)$);

\node[color=white] (Anpaux) at (0,-16) {$\mathbf{Pth}_{\boldsymbol{\mathcal{A}}^{(n)}}/{\mathrm{Ker}(\mathrm{CH}^{(n)})}$};
\node[color=white] (Anptaux) at (0,-18) {$\mathbf{Pth}_{\boldsymbol{\mathcal{A}}^{(n)}}/{\mathrm{Ker}(\mathrm{CH}^{(n)})}$};
\node [encercla, draw=blue!65!black,  fit=(A2paux)(A2ptaux)] {};

\node[color=blue!45!black] (Anp) at (0,-16) {$\llbracket \mathbf{Pth}_{\boldsymbol{\mathcal{A}}^{(n)}}\rrbracket$};
\node[color=blue!45!black] (Anpt) at (0,-18) {$\llbracket \mathbf{PT}_{\boldsymbol{\mathcal{A}}^{(n)}}\rrbracket$};

\node[rotate=90, color=blue!45!black] (Aniso) at (0,-17) {$\cong$};
\draw[ACliment, color=blue!45!black]  ($(Aniso)+(-.3,.7)$) to node [sloped, rotate=180, above] {
$\textstyle \mathrm{CH}^{\llbracket n \rrbracket}$
} ($(Aniso)+(-.3,-.7)$);
\draw[ACliment, color=blue!45!black]  ($(Aniso)+(.3,-.7)$) to node [sloped, below] {
$\textstyle \mathrm{ip}^{(\llbracket n \rrbracket, X)@}$
} ($(Aniso)+(.3,.7)$);

\node[] (A5aux) at (0,-19.5) {};
\draw[ACliment, color=blue!45!black]  ($(A5aux)+(0,1)$) to node [sloped, rotate=180, midway, fill=white] {
$\textstyle \mathrm{ip}^{(\llbracket n+1\rrbracket,\llbracket n\rrbracket)\sharp}$
} ($(A5aux)+(0,-1)$);
\draw[ACliment, color=blue!45!black]  ($(A5aux)+(.3,-1)$) to node [sloped, below, fill=white] {
$\textstyle \mathrm{tg}^{(\llbracket n\rrbracket,\llbracket n+1\rrbracket)}$
} ($(A5aux)+(.3,1)$);
\draw[ACliment, color=blue!45!black]  ($(A5aux)+(-.3,-1)$) to node [sloped, above, fill=white] {
$\textstyle \mathrm{sc}^{(\llbracket n\rrbracket, \llbracket n+1\rrbracket)}$
} ($(A5aux)+(-.3,1)$);

\node[] (ca) at (0,-20.85) {$\vdots$};

\node [encercla, draw=blue!45!black,  fit=(Anpaux)(Anptaux)] {};



\node[color=white] (B0aux) at (7,0){$\mathbf{Pth}_{\boldsymbol{\mathcal{A}}^{(n)}}/{\mathrm{Ker}(\mathrm{CH}^{(n)})}$};
\node [encercla, draw=red!95!black,  fit=(B0aux)] {} ;
\node[color=red!95!black] (B0) at (7,0) {$\mathbf{T}_{\Lambda}(Y)$};

\node[] (B0aux) at (7,-1.5) {};
\draw[ACliment, color=red!95!black]  ($(B0aux)+(0,1)$) to node [sloped, rotate=180, midway, fill=white] {
$\textstyle \mathrm{ip}^{([1],0)\sharp}$
} ($(B0aux)+(0,-1)$);
\draw[ACliment, color=red!95!black]  ($(B0aux)+(.3,-1)$) to node [sloped, below, fill=white] {
$\textstyle \mathrm{tg}^{(0,[1])}$
} ($(B0aux)+(.3,1)$);
\draw[ACliment, color=red!95!black]  ($(B0aux)+(-.3,-1)$) to node [sloped, above, fill=white] {
$\textstyle \mathrm{sc}^{(0,[1])}$
} ($(B0aux)+(-.3,1)$);

\node[color=white] (B1paux) at (7,-3) {$\mathbf{Pth}_{\boldsymbol{\mathcal{A}}^{(n)}}/{\mathrm{Ker}(\mathrm{CH}^{(n)})}$};
\node[color=white] (B1ptaux) at (7,-5) {$\mathbf{Pth}_{\boldsymbol{\mathcal{A}}^{(n)}}/{\mathrm{Ker}(\mathrm{CH}^{(n)})}$};

\node [encercla, draw=red!80!black,  fit=(B1paux)(B1ptaux)] {};

\node[color=red!80!black] (B1p) at (7,-3) {$[\mathbf{Pth}_{\boldsymbol{\mathcal{B}}}]$};
\node[color=red!80!black] (B1pt) at (7,-5) {$[\mathbf{PT}_{\boldsymbol{\mathcal{B}}}]$};

\node[rotate=90, color=red!80!black] (B1iso) at (7,-4) {$\cong$};
\draw[ACliment, color=red!80!black]  ($(B1iso)+(-.3,.7)$) to node [sloped, rotate=180, above] {
$\textstyle \mathrm{CH}^{[1]}$
} ($(B1iso)+(-.3,-.7)$);
\draw[ACliment, color=red!80!black]  ($(B1iso)+(.3,-.7)$) to node [sloped, below] {
$\textstyle \mathrm{ip}^{([1],Y)@}$
} ($(B1iso)+(.3,.7)$);

\node[] (B1aux) at (7,-6.5) {};
\draw[ACliment, color=red!80!black]  ($(B1aux)+(0,1)$) to node [sloped, rotate=180, midway, fill=white] {
$\textstyle \mathrm{ip}^{(\llbracket 2\rrbracket,[1])\sharp}$
} ($(B1aux)+(0,-1)$);
\draw[ACliment, color=red!80!black]  ($(B1aux)+(.3,-1)$) to node [sloped, below, fill=white] {
$\textstyle \mathrm{tg}^{([1],\llbracket 2\rrbracket)}$
} ($(B1aux)+(.3,1)$);
\draw[ACliment, color=red!80!black]  ($(B1aux)+(-.3,-1)$) to node [sloped, above, fill=white] {
$\textstyle \mathrm{sc}^{([1],\llbracket 2\rrbracket)}$
} ($(B1aux)+(-.3,1)$);

\node[color=white] (B2paux) at (7,-8){$\mathbf{Pth}_{\boldsymbol{\mathcal{A}}^{(n)}}/{\mathrm{Ker}(\mathrm{CH}^{(n)})}$};
\node[color=white] (B2ptaux) at (7,-10) {$\mathbf{Pth}_{\boldsymbol{\mathcal{A}}^{(n)}}/{\mathrm{Ker}(\mathrm{CH}^{(n)})}$};
\node [encercla, draw=red!65!black,  fit=(B2paux)(B2ptaux)] {};

\node[color=red!65!black] (B2p) at (7,-8) {$\llbracket \mathbf{Pth}_{\boldsymbol{\mathcal{B}}^{(2)}}\rrbracket$};
\node[color=red!65!black] (B2pt) at (7,-10) {$\llbracket \mathbf{PT}_{\boldsymbol{\mathcal{B}}^{(2)}}\rrbracket$};

\node[rotate=90, color=red!65!black] (B2iso) at (7,-9) {$\cong$};
\draw[ACliment, color=red!65!black]  ($(B2iso)+(-.3,.7)$) to node [sloped, rotate=180, above] {
$\textstyle \mathrm{CH}^{\llbracket 2 \rrbracket}$
} ($(B2iso)+(-.3,-.7)$);
\draw[ACliment, color=red!65!black]  ($(B2iso)+(.3,-.7)$) to node [sloped, below] {
$\textstyle \mathrm{ip}^{(\llbracket 2 \rrbracket, Y)@}$
} ($(B2iso)+(.3,.7)$);

\node[] (B2aux) at (7,-11.5) {};
\draw[ACliment, color=red!65!black]  ($(B2aux)+(0,1)$) to node [sloped, rotate=180, midway, fill=white] {
$\textstyle \mathrm{ip}^{(\llbracket 3\rrbracket, \llbracket 2\rrbracket)\sharp}$
} ($(B2aux)+(0,-1)$);
\draw[ACliment, color=red!65!black]  ($(B2aux)+(.3,-1)$) to node [sloped, below, fill=white] {
$\textstyle \mathrm{tg}^{(\llbracket 2\rrbracket,\llbracket 3\rrbracket)}$
} ($(B2aux)+(.3,1)$);
\draw[ACliment, color=red!65!black]  ($(B2aux)+(-.3,-1)$) to node [sloped, above, fill=white] {
$\textstyle \mathrm{sc}^{(\llbracket 2\rrbracket, \llbracket 3\rrbracket)}$
} ($(B2aux)+(-.3,1)$);

\node[] (ca) at (7,-12.85) {$\vdots$};

\node[] (B3aux) at (7,-14.5) {};
\draw[ACliment, color=red!50!black]  ($(B3aux)+(0,1)$) to node [sloped, rotate=180, midway, fill=white] {
$\textstyle \mathrm{ip}^{(\llbracket n\rrbracket, \llbracket n-1\rrbracket)\sharp}$
} ($(B3aux)+(0,-1)$);
\draw[ACliment, color=red!50!black]  ($(B3aux)+(.3,-1)$) to node [sloped, below, fill=white] {
$\textstyle \mathrm{tg}^{(\llbracket n-1\rrbracket, \llbracket n\rrbracket)}$
} ($(B3aux)+(.3,1)$);
\draw[ACliment, color=red!50!black]  ($(B3aux)+(-.3,-1)$) to node [sloped, above, fill=white] {
$\textstyle \mathrm{sc}^{(\llbracket n-1\rrbracket, \llbracket n\rrbracket)}$
} ($(B3aux)+(-.3,1)$);

\node[color=white] (Bnpaux) at (7,-16) {$\mathbf{Pth}_{\boldsymbol{\mathcal{A}}^{(n)}}/{\mathrm{Ker}(\mathrm{CH}^{(n)})}$};
\node[color=white] (Bnptaux) at (7,-18) {$\mathbf{Pth}_{\boldsymbol{\mathcal{A}}^{(n)}}/{\mathrm{Ker}(\mathrm{CH}^{(n)})}$};

\node[color=red!45!black] (Bnp) at (7,-16) {$\llbracket \mathbf{Pth}_{\boldsymbol{\mathcal{B}}^{(n)}}\rrbracket$};
\node[color=red!45!black] (Bnpt) at (7,-18) {$\llbracket \mathbf{PT}_{\boldsymbol{\mathcal{B}}^{(n)}} \rrbracket$};

\node[rotate=90, color=red!45!black] (Bniso) at (7,-17) {$\cong$};
\draw[ACliment, color=red!45!black]  ($(Bniso)+(-.3,.7)$) to node [sloped, rotate=180, above] {
$\textstyle \mathrm{CH}^{\llbracket n \rrbracket}$
} ($(Bniso)+(-.3,-.7)$);
\draw[ACliment, color=red!45!black]  ($(Bniso)+(.3,-.7)$) to node [sloped, below] {
$\textstyle \mathrm{ip}^{(\llbracket n \rrbracket, Y)@}$
} ($(Bniso)+(.3,.7)$);

\node[] (B5aux) at (7,-19.5) {};
\draw[ACliment, color=red!45!black]  ($(B5aux)+(0,1)$) to node [sloped, rotate=180, midway, fill=white] {
$\textstyle \mathrm{ip}^{(\llbracket n+1\rrbracket, \llbracket n\rrbracket)\sharp}$
} ($(B5aux)+(0,-1)$);
\draw[ACliment, color=red!45!black]  ($(B5aux)+(.3,-1)$) to node [sloped, below, fill=white] {
$\textstyle \mathrm{tg}^{(\llbracket n\rrbracket, \llbracket n+1\rrbracket)}$
} ($(B5aux)+(.3,1)$);
\draw[ACliment, color=red!45!black]  ($(B5aux)+(-.3,-1)$) to node [sloped, above, fill=white] {
$\textstyle \mathrm{sc}^{(\llbracket n\rrbracket, \llbracket n+1\rrbracket)}$
} ($(B5aux)+(-.3,1)$);

\node[] (ca) at (7,-20.85) {$\vdots$};

\node [encercla, draw=red!45!black,  fit=(Bnpaux)(Bnptaux)] {};


\node[] (M0) at (3.5,0) {};
\draw[ACliment]  ($(M0)+(-1,0)$) to node [above] 
{$f^{(0)\sharp}$} ($(M0)+(1,0)$);

\node[] (M1) at (3.5,-4) {};
\draw[ACliment]  ($(M1)+(-1,0)$) to node [above] 
{$f^{[1]@}$}  ($(M1)+(1,0)$);

\node[] (M2) at (3.5,-9) {};
\draw[ACliment]  ($(M2)+(-1,0)$) to node [above] 
{$f^{\llbracket 2\rrbracket@}$}  ($(M2)+(1,0)$);

\node[] (M2aux) at (3.5,-12.85) {$\cdots$};

\node[] (M3) at (3.5,-17) {};
\draw[ACliment]  ($(M3)+(-1,0)$) to node [above] 
{$f^{\llbracket n\rrbracket@}$}  ($(M3)+(1,0)$);

\node[] (M3aux) at (3.5,-20.85) {$\cdots$};

\end{tikzpicture}
\end{center}
\caption{A morphism of towers with the underlying Curry-Howard isomorphisms.}
\label{FA3}
\end{figure}

This finally, after verifying the compatibility with the structural morphisms of the towers, yields the morphism $\mathbf{f}^{\llbracket \bigcdot \rrbracket}$ from $\mathbb{A}$ to $\mathbb{B}$ (metaphorically speaking, $f^{(\bigcdot)}$ is the seed the morphism $\mathbf{f}^{\llbracket \bigcdot \rrbracket}$ grows up from) which will be graphically represented, in full, as shown in the diagrams of Figures~\ref{FA3} and~\ref{FA4}.
\begin{figure}
\begin{adjustbox}{width=\textwidth}
\begin{tikzpicture}
[ACliment/.style={-{To [angle'=45, length=5.75pt, width=4pt, round]}},
RHACliment/.style={right hook-{To [angle'=45, length=5.75pt, width=4pt, round]}, font=\scriptsize}, 
scale=.9
]

\tikzset{encercla/.style={draw=black, line width=.5pt, inner sep=0pt, rectangle, rounded corners}};

\node[] (A0g) at (0,3) {$X$};
\node[] (A0f) at (0,0) {\textcolor{white}{$\mathbf{T}_{\boldsymbol{\mathcal{E}}^{\boldsymbol{\mathcal{A}}^{(n)}}}(\mathbf{Pth}_{\boldsymbol{\mathcal{A}}^{(n)}})$}};
\node[] (A0) at (0,0) {$\mathbf{T}_{\Sigma}(X)$};

\node[] (A1g) at (4.5,3) {$\mathcal{A}$};
\node[] (A1f) at (4.5,0) {\textcolor{white}{$\mathbf{T}_{\boldsymbol{\mathcal{E}}^{\boldsymbol{\mathcal{A}}^{(n)}}}(\mathbf{Pth}_{\boldsymbol{\mathcal{A}}^{(n)}})$}};
\node[] (A1) at (4.5,0) {$\mathbf{T}_{\boldsymbol{\mathcal{E}}^{\boldsymbol{\mathcal{A}}}}(\mathbf{Pth}_{\boldsymbol{\mathcal{A}}})$};

\node[] (A2g) at (10.5,3) {$\mathcal{A}^{(n-1)}$};
\node[] (A2f) at (10.5,0) {\textcolor{white}{$\mathbf{T}_{\boldsymbol{\mathcal{E}}^{\boldsymbol{\mathcal{A}}^{(n-1)}}}(\mathbf{Pth}_{\boldsymbol{\mathcal{A}}^{(n-1)}})$}};
\node[] (A2) at (10.5,0) {$\mathbf{T}_{\boldsymbol{\mathcal{E}}^{\boldsymbol{\mathcal{A}}^{(n-1)}}}(\mathbf{Pth}_{\boldsymbol{\mathcal{A}}^{(n-1)}})$};

\node[] (Ang) at (16,3) {$\mathcal{A}^{(n)}$};
\node[] (An) at (16,0) {$\mathbf{T}_{\boldsymbol{\mathcal{E}}^{\boldsymbol{\mathcal{A}}^{(n)}}}(\mathbf{Pth}_{\boldsymbol{\mathcal{A}}^{(n)}})$};

\node[] (B0g) at (-2,1) {$\mathbf{T}_{\Lambda}(Y)$};
\node[] (B0f) at (-2,-2) {\textcolor{white}{$\mathbf{T}_{\boldsymbol{\mathcal{E}}^{\boldsymbol{\mathcal{B}}^{(n)}}}(\mathbf{Pth}_{\boldsymbol{\mathcal{B}}^{(n)}})$}};
\node[] (B0) at (-2,-2) {$\mathbf{T}_{\Lambda}(Y)_{\varphi}$};

\node[] (B1g) at (2.5,1) {$\mathrm{Pth}_{\mathcal{B},\varphi}$};
\node[] (B1f) at (2.5,-2) {\textcolor{white}{$\mathbf{T}_{\boldsymbol{\mathcal{E}}^{\boldsymbol{\mathcal{B}}^{(n)}}}(\mathbf{Pth}_{\boldsymbol{\mathcal{B}}^{(n)}})$}};
\node[] (B1) at (2.5,-2) {$\mathbf{T}^{\mathbf{f}^{(1)}}_{\boldsymbol{\mathcal{E}}^{\boldsymbol{\mathcal{B}}}}(\mathbf{Pth}_{\boldsymbol{\mathcal{B}}})$};

\node[] (B2g) at (8.5,1) {$\mathrm{Pth}_{\mathcal{B}^{(n-1)},\varphi}$};
\node[] (B2f) at (8.5,-2) {\textcolor{white}{$\mathbf{T}_{\boldsymbol{\mathcal{E}}^{\boldsymbol{\mathcal{B}}^{(n)}}}(\mathbf{Pth}_{\boldsymbol{\mathcal{B}}^{(n)}})$}};
\node[] (B2) at (8.5,-2) {$\mathbf{T}^{\mathbf{f}^{(n-1)}}_{\boldsymbol{\mathcal{E}}^{\boldsymbol{\mathcal{B}}^{(n-1)}}}(\mathbf{Pth}_{\boldsymbol{\mathcal{B}}^{(n-1)}})$};

\node[] (Bng) at (14,1) {$\mathrm{Pth}_{\mathcal{B}^{(n)},\varphi}$};
\node[] (Bn) at (14,-2) {$\mathbf{T}^{\mathbf{f}^{(n)}}_{\boldsymbol{\mathcal{E}}^{\boldsymbol{\mathcal{B}}^{(n)}}}(\mathbf{Pth}_{\boldsymbol{\mathcal{B}}^{(n)}})$};


\draw[ACliment]  (A0g) to node [left, pos=0.3] {$f^{(0)}$} (B0g);
\draw[ACliment]  (A1g) to node [left,  pos=0.3] {$f^{(1)}$} (B1g);
\draw[ACliment]  (A2g) to node [left,  pos=0.3] {$f^{(n-1)}$} (B2g);
\draw[ACliment]  (Ang) to node [left,  pos=0.3] {$f^{(n)}$} (Bng);

\draw[ACliment]  (B0g) to node [right] {
} ($(B0f)+(0,.4)$);

\draw[ACliment]  ($(B1f)+(-1.4,.15)$) to node [above, pos=.30 ] {
} ($(B0f)+(1,.15)$);
\draw[ACliment]  ($(B1f)+(-1.4,-.15)$) to node [below, pos=.30 ] {
} ($(B0f)+(1,-.15)$);
\draw[ACliment]  ($(B0f)+(1,0)$) to node [below, pos=.30 ] {
} ($(B1f)+(-1.4,0)$);
\draw[ACliment, rotate=30, gray!30!white]  ($(B1g)+(-.5,.15)$) to node [above, pos=.40 ] {
} ($(B0f)+(1.2,.15)$);
\draw[ACliment, rotate=30, gray!30!white]  ($(B1g)+(-.5,-.15)$) to node  [below, pos=.25 ] {
} ($(B0f)+(1.2,-.15)$);
\draw[ACliment]  (B1g) to node [right] {
} ($(B1f)+(0,.4)$);

\draw[ACliment]  ($(Bn)+(-1.6,.15)$) to node [above, pos=.30 ] {
} ($(B2f)+(1.9,.15)$);
\draw[ACliment]  ($(Bn)+(-1.6,-.15)$) to node [below, pos=.30 ] {
} ($(B2f)+(1.9,-.15)$);
\draw[ACliment]  ($(B2f)+(1.9,0)$) to node [below, pos=.30 ] {
} ($(Bn)+(-1.6,0)$);
\draw[ACliment, rotate=30, gray!30!white]  ($(Bng)+(-.5,.15)$) to node [above, pos=.40 ] {
} ($(B2f)+(1.2,.15)$);
\draw[ACliment, rotate=30, gray!30!white]  ($(Bng)+(-.5,-.15)$) to node  [below, pos=.25 ] {
} ($(B2f)+(1.2,-.15)$);
\draw[ACliment]  (B2g) to node [right] {
} ($(B2f)+(0,.4)$);

\draw[ACliment]  (Bng) to node [right] {
} ($(Bn)+(0,.4)$);


\draw[ACliment]  (A0g) to node [right] {
} ($(A0f)+(0,.4)$);

\draw[ACliment]  ($(A1f)+(-1.4,.15)$) to node [above, pos=.30 ] {
} ($(A0f)+(1,.15)$);
\draw[ACliment]  ($(A1f)+(-1.4,-.15)$) to node [below, pos=.30 ] {
} ($(A0f)+(1,-.15)$);
\draw[ACliment]  ($(A0f)+(1,0)$) to node [below, pos=.30 ] {
} ($(A1f)+(-1.4,0)$);
\draw[ACliment, rotate=30, gray!30!white]  ($(A1g)+(-.5,.15)$) to node [above, pos=.40 ] {
} ($(A0f)+(1.2,.15)$);
\draw[ACliment, rotate=30, gray!30!white]  ($(A1g)+(-.5,-.15)$) to node  [below, pos=.25 ] {
} ($(A0f)+(1.2,-.15)$);
\draw[ACliment]  (A1g) to node [right] {
} ($(A1f)+(0,.4)$);

\draw[ACliment]  ($(An)+(-1.6,.15)$) to node [above, pos=.30 ] {
} ($(A2f)+(1.9,.15)$);
\draw[ACliment]  ($(An)+(-1.6,-.15)$) to node [below, pos=.30 ] {
} ($(A2f)+(1.9,-.15)$);
\draw[ACliment]  ($(A2f)+(1.9,0)$) to node [below, pos=.30 ] {
} ($(An)+(-1.6,0)$);
\draw[ACliment, rotate=30, gray!30!white]  ($(Ang)+(-.5,.15)$) to node [above, pos=.40 ] {
} ($(A2f)+(1.2,.15)$);
\draw[ACliment, rotate=30, gray!30!white]  ($(Ang)+(-.5,-.15)$) to node  [below, pos=.25 ] {
} ($(A2f)+(1.2,-.15)$);
\draw[ACliment]  (A2g) to node [right] {
} ($(A2f)+(0,.4)$);

\draw[ACliment]  (Ang) to node [right] {
} ($(An)+(0,.4)$);

\node[] () at (6.75,0) {$\cdots$};
\node[] () at (19,0) {$\cdots$};

\draw[ACliment]  (A0f) to node [right, pos=0.4] {$f^{(0)\sharp}$} (B0f);
\draw[ACliment]  (A1f) to node [right,  pos=0.4] {$f^{[1]@}$} (B1f);
\draw[ACliment]  (A2f) to node [right,  pos=0.4] {$f^{\llbracket n-1\rrbracket@}$} (B2f);
\draw[ACliment]  (An) to node [right,  pos=0.4] {$f^{\llbracket n\rrbracket@}$} (Bn);

\draw[ACliment]  (A2g) to node [right] {
} ($(A2f)+(0,.4)$);

\end{tikzpicture}
\end{adjustbox}
\caption{A morphism from the tower $\mathbb{A}$ to the tower $\mathbb{B}$.}
\label{FA4}
\end{figure}

\begin{figure}[h]
\begin{center}
\begin{tikzpicture}
[ACliment/.style={-{To [angle'=45, length=5.75pt, width=4pt, round]}},
RHACliment/.style={right hook-{To [angle'=45, length=5.75pt, width=4pt, round]}, font=\scriptsize}, 
scale=.85
]

\node[] (T) at (0,0) {$\mathsf{Tw}$};
\node[] (LRws) at (-4,-2) {$\varprojlim \left(\mathsf{Rws}_{\mathfrak{d}}^{\llbracket \bigcdot \rrbracket}, F^{\llbracket \bigcdot^{+} \rrbracket}\right)$};
\node[] (LInt) at (4,-2) {$\varprojlim \left(
\int^{\mathsf{Sig}_{\mathfrak{d}}}(\bigcdot)\mathrm{CA},
\int^{\mathrm{Id}_{\mathsf{Sig}_{\mathfrak{d}}}}\eta^{(\bigcdot^{+})}
\right)$};
\node[] (Sig) at (0,-10) {$\mathsf{Sig}_{\mathfrak{d}}$};

\node[] (Rwsn) at (-1.5,-5) {
$\mathsf{Rws}_{\mathfrak{d}}^{\llbracket n \rrbracket}$};
\node[] (Rwsn1) at (-6.5,-7) {
$\mathsf{Rws}_{\mathfrak{d}}^{\llbracket n+1 \rrbracket}$};

\node[] (Intn) at (6.5,-5) {
$\int^{\mathsf{Sig}_{\mathfrak{d}}}\mathrm{nCA}$};
\node[] (Intn1) at (1.5,-7) {
$\int^{\mathsf{Sig}_{\mathfrak{d}}}\mathrm{(n+1)CA}$};

\draw[ACliment]  (LRws) to node [above] {
$\varprojlim E^{(\bigcdot)}$
} (LInt);

\draw[ACliment]  (Rwsn1) to node [midway, fill=white] {
$ F^{\llbracket n+1 \rrbracket}$
} (Rwsn);
\draw[ACliment]  (Intn1) to node [midway, fill=white] {
$\int^{\mathrm{Id}_{\mathsf{Sig}_{\mathfrak{d}}}}\eta^{(n+1)}$
} (Intn);

\draw[ACliment]  (Rwsn1) to node [pos=.4, below] {
$ E^{(n+1)}$
} (Intn1);
\draw[ACliment]  (Rwsn) to node [pos=.4, above] {
$ E^{(n)}$
} (Intn);

\draw[ACliment]  (T) to node [] {} (LInt);
\draw[ACliment]  (T) to node [] {} (LRws);
\draw[ACliment]  (LRws) to node [] {} (Rwsn);
\draw[ACliment]  (LRws) to node [] {} (Rwsn1);
\draw[ACliment]  (LInt) to node [] {} (Intn);
\draw[ACliment]  (LInt) to node [] {} (Intn1);

\draw[ACliment]  (Rwsn) to node [midway, fill=white] {$S^{(n)}_{\mathfrak{d}}$} (Sig);
\draw[ACliment]  (Rwsn1) to node [below left] {$S^{(n+1)}_{\mathfrak{d}}$} (Sig);

\draw[ACliment]  (Intn1) to node [midway, fill=white] {$\pi_{(n+1)\mathrm{CA}}$} (Sig);
\draw[ACliment]  (Intn) to node [below right] {$\pi_{n\mathrm{CA}}$} (Sig);

\end{tikzpicture}
\caption{From $\mathsf{Tw}$ to $\mathsf{Sig}_{\mathfrak{d}}$.}
\label{FAN}
\end{center}
\end{figure}

We let $\mathsf{Rws}_{\mathfrak{d}}^{\llbracket \bigcdot \rrbracket}$ stand for the family of categories 
$(\mathsf{Rws}_{\mathfrak{d}}^{\llbracket n \rrbracket})_{n\in\mathbb{N}}$. Then, since, for every $n\in \mathbb{N}$, there exists a forgetful functor $F^{\llbracket n+1 \rrbracket}$ from $\mathsf{Rws}_{\mathfrak{d}}^{\llbracket n+1 \rrbracket}$ to $\mathsf{Rws}_{\mathfrak{d}}^{\llbracket n \rrbracket}$, we obtain the projective system 
$(\mathsf{Rws}_{\mathfrak{d}}^{\llbracket\bigcdot\rrbracket},F^{\llbracket\bigcdot^{+}\rrbracket})$ in $\mathsf{Cat}$ and its projective limit $\varprojlim(\mathsf{Rws}_{\mathfrak{d}}^{\llbracket\bigcdot\rrbracket},F^{\llbracket\bigcdot^{+}\rrbracket})$. 

On the other hand, for every $n\in\mathbb{N}$, we have a contravariant functor 
$\mathrm{nCA}$ from $\mathsf{Sig}_{\mathfrak{d}}$ to $\mathsf{Cat}$, that assigns to a many-sorted signature $(S,\Sigma)$ the category $\mathsf{nCat}^{S}\mathsf{Alg}(\Sigma)$ of $S$-sorted $n$-categorial $\Sigma$-algebras. We let $\mathrm{(\bigcdot)CA}$ stand for $(\mathrm{nCA})_{n\in\mathbb{N}}$. Then, since, for every $n\in \mathbb{N}$, there exists a morphism $(\mathrm{Id}_{\mathsf{Sig}_{\mathfrak{d}}},\eta ^{(n+1)})$ from the split indexed category $(\mathsf{Sig}_{\mathfrak{d}},\mathrm{(n+1)CA})$ to the split indexed category $(\mathsf{Sig}_{\mathfrak{d}},\mathrm{nCA})$, we obtain the projective system $((\mathsf{Sig}_{\mathfrak{d}},\mathrm{(\bigcdot)CA}),(\mathrm{Id}_{\mathsf{Sig}_{\mathfrak{d}}},\eta ^{(\bigcdot^{+})}))$ in the category 
$\mathsf{SICat}$, of split indexed categories, and its projective limit 
$\varprojlim ((\mathsf{Sig}_{\mathfrak{d}},\mathrm{(\bigcdot)CA}),(\mathrm{Id}_{\mathsf{Sig}_{\mathfrak{d}}},\eta ^{(\bigcdot^{+})}))$. Now, since the Grothendieck construction, denoted by $\mathrm{Gr}$, is an equivalence from $\mathsf{SICat}$ to the category $\mathsf{SFib}$, of split fibrations, we obtain an equivalence $\mathrm{Gr}^{\boldsymbol{\omega}^{\mathrm{op}}}$ from $\mathsf{SICat}^{\boldsymbol{\omega}^{\mathrm{op}}}$ to $\mathsf{SFib}^{\boldsymbol{\omega}^{\mathrm{op}}}$. Hence, from the projective system 
$((\mathsf{Sig}_{\mathfrak{d}},\mathrm{(\bigcdot)CA}),(\mathrm{Id}_{\mathsf{Sig}_{\mathfrak{d}}},\eta ^{(\bigcdot^{+})}))$ in $\mathsf{SICat}$, we obtain the projective system 
$((\int^{\mathsf{Sig}_{\mathfrak{d}}}(\bigcdot)\mathrm{CA},\pi_{(\bigcdot)\mathrm{CA}}),(\int^{\mathrm{Id}_{\mathsf{Sig}_{\mathfrak{d}}}}\eta^{(\bigcdot^{+})},\mathrm{Id}_{\mathsf{Sig}_{\mathfrak{d}}}))$ in the category $\mathsf{SFib}$ and its projective limit 
$\varprojlim ((\int^{\mathsf{Sig}_{\mathfrak{d}}}(\bigcdot)\mathrm{CA},\pi_{(\bigcdot)\mathrm{CA}}),(\int^{\mathrm{Id}_{\mathsf{Sig}_{\mathfrak{d}}}}\eta^{(\bigcdot^{+})},\mathrm{Id}_{\mathsf{Sig}_{\mathfrak{d}}}))$. Moreover, for every $n\in\mathbb{N}$, there exists a functor $S^{\llbracket n \rrbracket}_{\mathfrak{d}}$ from $\mathsf{Rws}_{\mathfrak{d}}^{\llbracket n \rrbracket}$ to $\mathsf{Sig}_{\mathfrak{d}}$ that assigns to an $S$-sorted rewriting system $\boldsymbol{\mathcal{A}}^{(n)}$ of $n$-th order its underlying many-sorted  signature $(S,\Sigma^{\mathbcal{A}^{(n)}})$, and there exists an embedding $E^{(n)}$ from $\mathsf{Rws}_{\mathfrak{d}}^{\llbracket n \rrbracket}$ to $\int^{\mathsf{Sig}_{\mathfrak{d}}}\mathrm{nCA}$, that assigns to an $S$-sorted rewriting system $\boldsymbol{\mathcal{A}}^{(n)}$ of $n$-th order the ordered pair whose first coordinate is the underlying $S$-sorted rewriting system of zeroth order, i.e., $\boldsymbol{\mathcal{A}}^{(0)}=(S,\Sigma,X)$, and whose second coordinate is the free partial $\Sigma^{\boldsymbol{\mathcal{A}}^{(n)}}$-algebra $\mathbf{T}_{\boldsymbol{\mathcal{E}}^{(n)}}(\mathbf{Pth}_{\boldsymbol{\mathcal{A}}^{(n)}})$, which is an $S$-sorted $n$-categorial $\Sigma$-algebra. We let $E^{(\bigcdot)}$ stand for $(E^{(n)})_{n\in\mathbb{N}}$. Then $E^{(\bigcdot)}$ is a morphism from the projective system $(\mathsf{Rws}_{\mathfrak{d}}^{\llbracket \bigcdot \rrbracket},F^{\llbracket \bigcdot^{+} \rrbracket})$ to the projective system 
$((\int^{\mathsf{Sig}_{\mathfrak{d}}}(\bigcdot)\mathrm{CA},\pi_{(\bigcdot)\mathrm{CA}}),(\int^{\mathrm{Id}_{\mathsf{Sig}_{\mathfrak{d}}}}\eta^{(\bigcdot^{+})},\mathrm{Id}_{\mathsf{Sig}_{\mathfrak{d}}}))$ and we let 
$\varprojlim E^{(\bigcdot)}$ stand for the canonical morphism from the projective limit of the first projective system to the projective limit of the second projective system.
Then, after defining, for every tower $\mathbb{A}$ the identity at $\mathbb{A}$ and the composition of morphisms between towers, we obtain the category $\mathsf{Tw}$ of towers and prove that there exists a functor from $\mathsf{Tw}$ to $\varprojlim (\mathsf{Rws}_{\mathfrak{d}}^{\llbracket \bigcdot \rrbracket},F^{\llbracket \bigcdot^{+} \rrbracket})$, the properties of which will be investigated. We graphically represent the above results in the diagram of Figure~\ref{FAN}.

Our \textsf{third objective} will be to investigate the geometry of the towers and their classifying spaces through the investigation of the geometry of the $S$-sorted $\omega$-categorial $\Sigma$-algebras and their classifying spaces, but this matter will require a deeper understanding of the issues being discussed here and we will leave this for future versions of the work.

Let us remark that in the above definition of the morphisms between towers, we have chosen to take the derivors in the basis of the recursive definition. It is worth noting that other, more complicated, morphisms could be chosen (see~\cite{CVST10b}) in the basis of the recursion. 
Moreover, by defining a suitable notion of transformation between derivors one can equip the category $\mathsf{Sig}_{\mathfrak{d}}$ with a structure of $2$-category (this was, in fact, already done in~\cite{CVST10b} for polyderivors, also called morphisms of Fujiwara there, of which the derivors are particular cases).

\section{Related work}

In the discussion that follows, we examine the relationship of our work to the seminal contributions of Hotz and Benson and those of Street and Burroni. Our focus is limited to the quoted authors, as they are the pioneers in the specific field we are addressing. However, we will also mention other scholars such as Power, who follows Street's approach, and Giraud, Malbos and Mimram, who follow Burroni's approach.

\subsection{The work of Hotz and Benson}
 
In~\cite{Hotz65, Hotz65a}, whose English title would be ``Algebraisation of the synthesis problem of switching circuits $\mathrm{I,II}$'', which together constitute Hotz's Habilitation Thesis, the author addresses challenges inherent in the synthesis of switching circuits by providing them with an algebraic framework. To achieve this, Hotz introduces the theory of $X$-categories, which are particular $2$-categories and are currently called strict monoidal categories. Actually, from a set $\mathfrak{A}$, whose elements represent the primitive building blocks,  Hotz obtains a free $X$-category $\mathfrak{F}(\mathfrak{A})$. The elements of $\mathfrak{F}(\mathfrak{A})$ represent the synthesis of a circuit from the primitive building blocks and each synthesis of a circuit has a correspondence in $\mathfrak{F}(\mathfrak{A})$. (Let us add that, in~\cite{Hotz65}, Hotz also introduced planar nets, currently called string diagrams, to represent the morphisms in strict monoidal categories. In other words, Hotz establishes the equivalence between the theory of strict monoidal categories and that of planar nets, i.e., a coherence theorem). Moreover, Hotz defines free $D$-categories,  i.e., free categories with direct products as quotient categories, $\mathfrak{F}(\mathfrak{A})/\mathfrak{R}$, given by generators 
$\mathfrak{A}$ and relations $\mathfrak{R}$. 

In~\cite{Hotz66}, whose English title would be ``Unambiguity and ambiguity of formal languages'', Hotz uses the strict monoidal categories and string diagrams to define the derivations in semi-Thue systems ($\equiv$ word rewriting system) and Chomsky grammars. Moreover, he defines the notion of homomorphism between Chomsky grammars. To the best of our knowledge, we think that with Hotz begins the categorial treatment of derivations in a word rewriting system.

Our work shares similarities and differences with that of Hotz. We both propose the definition of free categorial algebras from rewriting systems (Hotz for word rewriting systems, we for many-sorted term rewriting systems of order $n$). In fact, Hotz focuses his attention on defining $X$-categories and $D$-categories given by generators and relations, while we address the more general problem of proving the existence of many-sorted $n$-categorial $\Sigma$-algebra, for a many-sorted signature $\Sigma$, from  generators and relations.

Hotz's papers~\cite{Hotz65, Hotz65a, Hotz66}, which came to our attention in the winter of 2023, after we had finished the current version, have had no direct influence on our work, but, after studying them, we may regard them as precursors of our own.


In~\cite{Ben70}, Benson defines a syntax as a category of strings and derivations between strings generated by a semi-Thue system, where a semi-Thue system is given by an alphabet, a set of productions and a starting axiom, i.e., a non-empty word over the alphabet. Moreover, Benson defines a derivation as a sequence of strings obtained by repeatedly applying productions in context. This definition coincides with the definition of first-order path that we propose in our work, when the latter is particularized to the single-sorted signature of monoids. He also defines two full subcategories, determined by the starting axiom and a set of terminal strings. This corresponds to the definition of our set of paths, to paths starting on a given term, and paths starting on a given term an ending on a set of terms, respectively.

In addition, in~\cite{Ben70}, Benson defines the semantic domain to be a category of sets and functions having cartesian products. Then, Benson defines the interpretations of a syntax as contravariant functors from the syntax to the semantic domain, generated by a correspondence between production of the semi-Thue system and certain functions. Following this, Benson proves that there exists a Galois connection between congruence relations and classes of interpretations. This illuminates the relationship between syntax and semantics. Benson devotes special attention to the smallest congruence of interest, similarity, which, by the way, corresponds to the class of all interpretations. Two derivations are considered similar if they can be transformed into each other by applying a series of switchings between the productions occurring on non-interacting subwords. Benson states that each similarity class contains a canonical derivation following a leftmost derivation strategy and that two derivations are similar if, and only if, they are identically interpreted by each possible interpretation. This notion of similarity coincides with the description of the kernel of the Curry-Howard mapping in our work. The remaining statements refer to the description of the normalized path in each class and to the universal property of the quotient of the syntax category by the similarity. According to Benson, similarity is the smallest congruence compatible with parallel exchange. In this respect, for us two derivations are similar if they have the same image under the Curry-Howard mapping.

The rest of the paper delves into other important issues. In particular, Benson investigates the quotients of a syntax by congruences containing the similarity congruence to understand the relationship between them and certain subclasses of interpretations.

Before describing the content of Benson's paper~\cite{Ben74}, let us point out that the translations between formal languages that Benson considers are the syntax directed translations and that for him translations must be effective to be of interest in current formal studies of natural and programming languages. In this respect he says that ``A natural method of obtaining effectiveness is to define it first in terms of the grammar generating the syntactic structures of the first language system and then by induction, \emph{of one form or another} [we emphasize], on the entire set of syntactic structures''. This approach coincides, in the particular case of the (single-sorted) signature of monoids, with the one we will follow when we will define, in the third part of our work, the translations between many-sorted term rewriting systems. 

In~\cite{Ben74}, Benson begins by presenting those category-theoretic notions wich will be used in the paper. Next, following~\cite{EW67}, he defines the notion of (single-sorted) algebraic theory and shows that they are $X$-categories. After this, Benson presents a generalization of the algebraic theories of Lawvere in~\cite{Law03}, which he calls $\mathbf{\Omega}^{\star}$-theories, and shows that they are also $X$-categories (let us point out that many-sorted algebras were introduced by Higgins in~\cite{Hig63} and  many-sorted algebraic theories by Bénabou in~\cite{jb68}). $\mathbf{\Omega}^{\star}$-theories are the algebraic theories of many-sorted algebras. Afterwards, Benson assigns to a rewriting system $(\Sigma,P)$ an $X$-category $\mathbf{F}_{(\Sigma,P)}$ whose objects are the words on $\Sigma$ and whose morphisms from a word to another are derivations in which every step of the process is completely 
specified, and notices that ``For linguistic and semantic reasons one concludes that there are too many derivations when completely specified, and certain equivalence classes of derivations for the proper notion of syntactic structure''. From this Benson concludes that ``$\mathbf{F}_{(\Sigma,P)}/\sim$, which is 
$\mathbf{F}_{(\Sigma,P)}$ factored by the similarity relationship $\sim$, is the proper category in which to study the syntactic structure imposed on $\Sigma^{\star}$ by $(\Sigma,P)$''. $\mathbf{F}_{(\Sigma,P)}/\sim$ is also an $X$-category. Following this, Benson formalizes the notion of algebraic interpretation from an $X$-category of derivations $\mathbf{X}$ to an $\mathbf{\Omega}^{\star}$-theory $\mathbf{T}$ as an $X$-functor from $\mathbf{X}$ to $\mathbf{T}$ and, on such a notion, he gives a definition of semantic-preserving translations between two language systems. Then Benson, by abstraction from the notion of the category of a generalized sequential machine, defines the notion of translations category, which are spans of $X$-categories, and shows that for derivation systems $\mathbf{X}_{j}$ generated by grammars
$(\Sigma_{j},V_{j},P_{j},\zeta_{j})$, $j=1,2$, algebraic interpretations $I_{j}\colon \mathbf{X}_{j}\mor \mathbf{T}$, $j=1,2$, in an $\mathbf{\Omega}^{\star}$-theory $\mathbf{T}$ for which $I_{1}(\zeta_{1})= I_{2}(\zeta_{2})$, if $\mathbf{X}_{0}$ is a translation category with $X$-functors $F_{j}\colon \mathbf{X}_{0}\mor \mathbf{X}_{j}$,  $j=1,2$, and $\eta$ is a natural transformation from $F_{2}\circ I_{2}$ to $F_{1}\circ I_{1}$, satisfying a certain condition, then the translations induced by the translation category $\mathbf{X}_{0}$ are semantic-preserving. Finally, Benson, as application of the result on semantic preservation by a natural transformation, shows that generalized sequential machine maps form translation categories and gives effective sufficient conditions for semantic preservation of such translation categories.



Finally, in~\cite{db75}, Benson introduces some basic algebraic structures within the categories of derivations determined by a rewriting system. More explicitly, Benson, in~\cite{db75}, begins by defining the category of derivations $\mathbf{D}$ for an indexed rewriting system $(A,P)$, where $A$ is a set and $P$ a mapping from a set $J$ to $A^{\star}\times A^{\star}$, the cartesian product of the underlying set of the free monoid $\mathbf{A}^{\star}$ on $A$ with itself. Following this he defines the relation of similarity between derivations which can be given in several equivalent ways. One is by interchanges (inessential exchanges). A second is by semantic considerations. A third is by congruences on $\mathbf{D}$, and this is the method used by Benson in his work. The congruential definition of similarity requires the consideration of a family of pairs of endofunctors, concatenation (of words) on the right and concatenation (of words) on the left, on the category $\mathbf{D}$, and is the least congruence on $\mathbf{D}$ containing a suitable set of pairs of derivations in the definition of which the aforementioned endofunctors occur. Then, after establishing a series of technical results, Benson states that $\mathbf{S}$, the quotient of $\mathbf{D}$ modulo similarity, which he calls the syntax category, is a free $X$-category, i.e., a free strict monoidal category, which is a nice example of a $2$-category. Benson shows, in particular, that syntax categories are very useful for investigating (word) rewriting systems.


Two main differences emerge when we compare our work with those of Hotz and Benson. Firstly, the species of  algebraic structure considered by Hotz in~\cite{Hotz66} and Benson in~\cite{db75} is that of free monoids, while ours is richer: free many-sorted $\Sigma$-algebras, for a many-sorted signature $\Sigma$. This greater generality will be reflected in the equations defining the generalization of the notion of category, where the operations from the original signature behave as functors for the composition of derivations. Secondly, the methods proposed by Benson are semantic-based, whilst ours are syntactic-based, assigning, in a functional way, to every path a term in an enlarged signature that encapsulates every interesting transformation process occurring in the original path. Actually, the kernel of the just mentioned mapping will be the central object of our study. Moreover, in contrast to Benson's work, this kernel will not be defined by a minimality condition but, simply, as the kernel of a totally defined mapping. In fact, the term associated to the original path will be used later to provide a canonical representative of each path class, where every transformation is performed according to the interpretation of the operations.

On the other hand, the similarities become evident in other directions. Firstly, the importance of taking quotients in the derivation process to overcome the limitations of a sequential definition, and, secondly, the algebraic structure enclosed in this quotient and its good behaviour with respect to the sequential operator. We think that these two properties are universal in any derivation process between structured terms. 



Benson's work, of which we became aware in the fall of 2018, is at the origin of ours, and, partially, provides the solid foundation on which we have built it. In particular, his work on the basic algebraic structures in categories of derivations was crucial and the trigger for our own work. We are deeply grateful to him for this. 
We would also like to express our gratitude to G\'{e}cseg and Steinby~\cite{GS84}, Schmidt~\cite{sch66,sch70}, Burmeister and Schmidt~\cite{bs67}, Burmeister~\cite{pb70, pb82, pb86, pb93, pb02} and Burmeister and Wojdy{\l}o~\cite{pbw92a, pbw92b, pbw96} for their just mentioned contributions, as we have relied on them in part for our work.





\subsection{The work of Street and Burroni}

Higher-order rewriting systems arise in the context of higher-order categories when the classical notion of presentation of a universal algebra by generators and relations is generalized, by using computads or, equivalently, polygraphs, up to that of presentation of an $n$-category. In~\cite{Str76}, after defining a graph $G = (G_{0},G_{1},d_{0},d_{1})$ as trim if $[d_{0},d_{1}]\colon G_{1}\coprod G_{1}\mor G_{0}$ is surjective, Street defines, in the $2$-dimensional case, a \emph{computad} $\mathcal{G}$ as consisting of a graph $\bb{\mathcal{G}}$ together with, for each ordered pair $A$, $B$ of objects of $\bb{\mathcal{G}}$, a trim graph $\mathcal{G}(A,B)$ such that $\mathcal{G}(A,B)_{0}$ is a subset of $(\mathcal{F}\bb{\mathcal{G}})(A,B)$, where $\mathcal{F}\bb{\mathcal{G}}$ is the free category on 
$\bb{\mathcal{G}}$. He shown that the category of $2$-categories is monadic over the category of computads and defines a presentation of a $2$-category $\mathcal{A}$ as a pair of computads 
$\mathcal{G}$, $\mathcal{H}$ together with a coequalizer $W \colon \mathcal{FG} \mor \mathcal{A}$ of 
$U, V \colon \mathcal{FH} \mor \mathcal{FG}$, where $\mathcal{FG}$, $\mathcal{FH}$ are the $2$-categories generated by $\mathcal{G}$, $\mathcal{H}$. This was later generalised by Power, in~\cite{Pow91}, to the 
$n$-dimensional case. In~\cite{Bur91, Bur93} and independently of Street and Power, Burroni defines an $n$-\emph{polygraph} as consisting of a diagram, shown in Figure~\ref{FA6},
\begin{figure}
\begin{center}
\begin{tikzpicture}
[ACliment/.style={-{To [angle'=45, length=5.75pt, width=4pt, round]}},
RHACliment/.style={right hook-{To [angle'=45, length=5.75pt, width=4pt, round]}, font=\scriptsize}, 
scale=.95
]

\node[] (A) at (0,0) {$E_{0}$};
\node[] (Ap) at (0,2) {$E_{0}$};
\draw[double equal sign distance]  (Ap) to node [right] {} (A);

\node[] (B) at (2,0) {$E^{\ast}_{1}$};
\node[] (Bp) at (2,2) {$E_{1}$};

\draw[ACliment]  ($(B)+(-.3,.15)$) to node [above, pos=.40 ] {
$\scriptstyle s^{\ast}_{0}$
} ($(A)+(.3,.15)$);
\draw[ACliment]  ($(B)+(-.3,-.15)$) to node [below, pos=.40 ] {
$\scriptstyle t^{\ast}_{0}$
} ($(A)+(.3,-.15)$);
\draw[ACliment, rotate=45]  ($(Bp)+(-.3,.15)$) to node [above, pos=.40, fill=white ] {
$\scriptstyle s_{0}$
} ($(A)+(.3,.15)$);
\draw[ACliment, rotate=45]  ($(Bp)+(-.3,-.15)$) to node  [below, pos=.25, fill=white ] {
$\scriptstyle t_{0}$
} ($(A)+(.3,-.15)$);
\draw[ACliment]  (Bp) to node [right, pos=.3] {
$\scriptstyle i_{1}$} (B);

\node[] (C) at (4,0) {$E^{\ast}_{2}$};
\node[] (Cp) at (4,2) {$E_{2}$};

\draw[ACliment]  ($(C)+(-.3,.15)$) to node [above, pos=.40 ] {
$\scriptstyle s^{\ast}_{1}$
} ($(B)+(.3,.15)$);
\draw[ACliment]  ($(C)+(-.3,-.15)$) to node [below, pos=.40 ] {
$\scriptstyle t^{\ast}_{1}$
} ($(B)+(.3,-.15)$);
\draw[ACliment, rotate=45]  ($(Cp)+(-.3,.15)$) to node [above, pos=.40,fill=white ] {
$\scriptstyle s_{1}$
} ($(B)+(.3,.15)$);
\draw[ACliment, rotate=45]  ($(Cp)+(-.3,-.15)$) to node  [below, pos=.25, fill=white ] {
$\scriptstyle t_{1}$
} ($(B)+(.3,-.15)$);
\draw[ACliment]  (Cp) to node [right, pos=.3] {
$\scriptstyle i_{2}$} (C);

\node[] () at (5,1) {$\cdots$};

\node[] (Dv) at (6,0) {$E^{\ast}_{n-2}\quad$};

\node[] (D) at (8,0) {$\quad E^{\ast}_{n-1}$};
\node[] (Dp) at (8,2) {$\quad E_{n-1}$};

\draw[ACliment]  ($(D)+(-.3,.15)$) to node [above, pos=.40 ] {
$\scriptstyle s^{\ast}_{n-2}$
} ($(Dv)+(.3,.15)$);
\draw[ACliment]  ($(D)+(-.3,-.15)$) to node [below, pos=.40 ] {
$\scriptstyle t^{\ast}_{n-2}$
} ($(Dv)+(.3,-.15)$);
\draw[ACliment, rotate=45]  ($(Dp)+(-.3,.15)$) to node [above, pos=.40, fill=white ] {
$\scriptstyle s_{n-2}$
} ($(Dv)+(.3,.15)$);
\draw[ACliment, rotate=45]  ($(Dp)+(-.3,-.15)$) to node  [below, pos=.25 , fill=white] {
$\scriptstyle t_{n-2}$
} ($(Dv)+(.3,-.15)$);
\draw[ACliment]  (Dp) to node [right, pos=.3] {
$\scriptstyle i_{n-1}$} (D);

\node[] (Ep) at (10,2) {$\quad E_{n}$};

\draw[ACliment, rotate=45]  ($(Ep)+(-.3,.15)$) to node [above, pos=.40 ,fill=white] {
$\scriptstyle s_{n-1}$
} ($(D)+(.3,.15)$);
\draw[ACliment, rotate=45]  ($(Ep)+(-.3,-.15)$) to node  [below, pos=.25 , fill=white] {
$\scriptstyle t_{n-1}$
} ($(D)+(.3,-.15)$);

\end{tikzpicture}
\caption{An $n$-polygraph.}
\label{FA6}
\end{center}
\end{figure}
such that, for each $k$ in $n$, it happens that 
\begin{itemize}
\item the diagram
\begin{tikzpicture}
[ACliment/.style={-{To [angle'=45, length=5.75pt, width=4pt, round]}},
RHACliment/.style={right hook-{To [angle'=45, length=5.75pt, width=4pt, round]}, font=\scriptsize}, 
scale=.85, baseline=-.9mm
]

\node[] (A) at (0,0) {$E_{0}$};
\node[] (B) at (1.5,0) {$E^{\ast}_{1}$};

\draw[ACliment]  ($(B)+(-.3,.15)$) to node [above, pos=.40 ] {
} ($(A)+(.3,.15)$);
\draw[ACliment]  ($(B)+(-.3,-.15)$) to node [below, pos=.40 ] {
} ($(A)+(.3,-.15)$);

\node[] (C) at (3,0) {$\cdots$};

\draw[ACliment]  ($(C)+(-.3,.15)$) to node [above, pos=.40 ] {
} ($(B)+(.3,.15)$);
\draw[ACliment]  ($(C)+(-.3,-.15)$) to node [below, pos=.40 ] {
} ($(B)+(.3,-.15)$);

\node[] (D) at (4.5,0) {$E^{\ast}_{k}$};
\draw[ACliment]  ($(D)+(-.3,.15)$) to node [above, pos=.40 ] {
} ($(C)+(.3,.15)$);
\draw[ACliment]  ($(D)+(-.3,-.15)$) to node [below, pos=.40 ] {
} ($(C)+(.3,-.15)$);
\end{tikzpicture}
is a $k$-category;
\item the diagram
\begin{tikzpicture}
[ACliment/.style={-{To [angle'=45, length=5.75pt, width=4pt, round]}},
RHACliment/.style={right hook-{To [angle'=45, length=5.75pt, width=4pt, round]}, font=\scriptsize}, 
scale=.85, baseline=-.9mm
]

\node[] (A) at (0,0) {$E_{0}$};
\node[] (B) at (1.5,0) {$E^{\ast}_{1}$};

\draw[ACliment]  ($(B)+(-.3,.15)$) to node [above, pos=.40 ] {
} ($(A)+(.3,.15)$);
\draw[ACliment]  ($(B)+(-.3,-.15)$) to node [below, pos=.40 ] {
} ($(A)+(.3,-.15)$);

\node[] (C) at (3,0) {$\cdots$};

\draw[ACliment]  ($(C)+(-.3,.15)$) to node [above, pos=.40 ] {
} ($(B)+(.3,.15)$);
\draw[ACliment]  ($(C)+(-.3,-.15)$) to node [below, pos=.40 ] {
} ($(B)+(.3,-.15)$);

\node[] (D) at (4.5,0) {$E^{\ast}_{k}$};
\draw[ACliment]  ($(D)+(-.3,.15)$) to node [above, pos=.40 ] {
} ($(C)+(.3,.15)$);
\draw[ACliment]  ($(D)+(-.3,-.15)$) to node [below, pos=.40 ] {
} ($(C)+(.3,-.15)$);

\node[] (E) at (6,0) {$\quad E_{k+1}$};
\draw[ACliment]  ($(E)+(-.3,.15)$) to node [above, pos=.40 ] {
} ($(D)+(.3,.15)$);
\draw[ACliment]  ($(E)+(-.3,-.15)$) to node [below, pos=.40 ] {
} ($(D)+(.3,-.15)$);
\end{tikzpicture}
is a $(k+1)$-graph.
\end{itemize}
Polygraphs were introduced by Burroni for studying higher dimensional word problems. In fact, Burroni, after saying that the word problem on a monoid has two natural generalizations to categories and universal algebras, and that in the first case ``the words become the `paths' and the equality problems take the form of diagram commutation problems'', while, in the second case, ``the words become the `terms' and the rewriting systems set the rules for their equality'', unifies both extensions by stating that ``the rewriting problem for terms is nothing but a $2$-dimensional path problem in a 2-category. This observation leads to the general problem for $n$-paths in an $n$-category, or even in an $\infty$-category. A lot of computations made by categoricians are $1$-, $2$- or $3$-dimensional computations, and in fact $n$-dimensional computations take place in an $(n + 1)$-category''. In short, Burroni's investigations amount to considering that the computations in $\infty$-categories are equivalent to constructions of homotopies between expressions.


In Burroni's work, the passage from one layer to the next is done by representing the paths of the previous strata as terms and then taking the appropriate quotients so that the categorial axioms are fulfilled.  Since the only structure involved is that of category, all terms are constructed as formal compositions according to the generators, their sources, and targets.


Readers interested in further exploring these topics may wish to read, in the $2$-dimensional case, \cite{jb67,sM98}, in the $n$-dimensional case, \cite{Ehr64,Lei01}, and in the $\infty$-dimensional case, \cite{Lan21,Lei04,RV16}. The compendium~\cite{abgmmm23}, comprehensively encompasses numerous results put forth by the school of polygraph theorists led by Burroni.

Note that there is an analogy between higher-order rewriting and algebraic topology. Specifically, if we think of terms as points, then derivations between terms will be paths. Continuing with the analogy, it will not be enough to know how to make transitions between points, but also between paths (i.e., $2$-dimensional morphisms) and between paths of paths (i.e., $3$-dimensional morphisms), and so on, until we build a rewriting system of the desired order~\cite{Mim14}. This is why Burroni's recursive definition relies on the polygraphic structure of the previous layers to define the next layer. Recall that these derivations will be discrete, and the edges of the $(n+1)$-polygraph will act as the generators of the set of $n$-paths.

The usefulness of higher-order rewriting systems has been highlighted in some recent work, e.g.,~\cite{CM15, Gui19, GM09, Laf03, Mim14}, in which the authors mainly focus their attention on the investigation of algebraic presentations in the sense of Lawvere~\cite{Law03}. Moreover, one of the constants in these investigations consists in obtaining convergent (i.e., terminating and confluent) rewriting systems, which makes it possible to establish a normal form of the terms modulo the congruence generated by the rewriting rules.

With regard to the just mentioned investigations, several challenges have to be faced. One of them has to do with the conceptual and notational complexity of the structures treated. Another is related to the model used to describe the objects of interest. We have observed that the classes of terms relative to a congruence of interest are always the most difficult, but necessary, choice to be made in order to iterate the process. Once the choice is made, geometrical models are considered to represent the cells, giving priority to diagrammatic representations of the $n$-cells or to wire diagrams (up to planar transformations), as in~\cite{Laf03, Mim14}. 

As for the problem of higher-order rewriting, there are similarities between the approach using the notion of a polygraph and our approach. Note the similarity of the representation of a second-order rewriting system in our approach, as shown in Figure~\ref{FA7}, and the representation of it by a $2$-polygraph. There are also similarities with respect to the recursive approach that uses generators at each layer. These generators (pairs of directed terms) help shape the notion of a path, ultimately the free category associated with the set of generators. And, as in Burroni's case, precautions will have to be taken in the definition of these generators, mainly conditions on the sources and targets so that the 
$(n+1)$-cells are defined on $n$-parallel cells. Another similarity has to do with the need to take quotients for everything to work as it should, since these paths, by themselves, do not satisfy the categorial equations. 

\begin{figure}
\begin{center}
\begin{tikzpicture}
[ACliment/.style={-{To [angle'=45, length=5.75pt, width=4pt, round]}},
RHACliment/.style={right hook-{To [angle'=45, length=5.75pt, width=4pt, round]}, font=\scriptsize}, 
scale=.92
]

\tikzset{encercla/.style={draw=black, line width=.5pt, inner sep=0pt, rectangle, rounded corners}};


\node[] (A) at (0,0) {$\mathbf{T}_{\Sigma}(X)$};

\node[] (Ap) at (0,2) {$X$};
\draw[ACliment]  (Ap) to node [right] {
} ($(A)+(0,.4)$);


\node[] (B) at (4,0) {$\mathbf{T}_{\boldsymbol{\mathcal{E}}^{\boldsymbol{\mathcal{A}}}}(\mathbf{Pth}_{\boldsymbol{\mathcal{A}}})$};
\node[] (Bp) at (4,2) {$\mathcal{A}$};

\draw[ACliment]  ($(B)+(-1.3,.15)$) to node [above, pos=.30 ] {
} ($(A)+(.8,.15)$);
\draw[ACliment]  ($(B)+(-1.3,-.15)$) to node [below, pos=.30 ] {
} ($(A)+(.8,-.15)$);
\draw[ACliment, rotate=45]  ($(Bp)+(-.3,.15)$) to node [above, pos=.40 ] {
} ($(A)+(.7,.15)$);
\draw[ACliment, rotate=45]  ($(Bp)+(-.3,-.15)$) to node  [below, pos=.25 ] {
} ($(A)+(.7,-.15)$);
\draw[ACliment]  (Bp) to node [right] {
} ($(B)+(0,.4)$);


\node[] (C) at (8,0) {$\mathbf{T}_{\boldsymbol{\mathcal{E}}^{\boldsymbol{\mathcal{A}}^{(2)}}}(\mathbf{Pth}_{\boldsymbol{\mathcal{A}}^{(2)}})$};
\node[] (Cp) at (8,2) {$\mathcal{A}^{(2)}$};

\draw[ACliment]  ($(C)+(-1.65,.15)$) to node [above, pos=.30 ] {
} ($(B)+(1.2,.15)$);
\draw[ACliment]  ($(C)+(-1.65,-.15)$) to node [below, pos=.30 ] {
} ($(B)+(1.2,-.15)$);
\draw[ACliment, rotate=45]  ($(Cp)+(-.3,.15)$) to node [above, pos=.40 ] {
} ($(B)+(.7,.15)$);
\draw[ACliment, rotate=45]  ($(Cp)+(-.3,-.15)$) to node  [below, pos=.25 ] {
} ($(B)+(.7,-.15)$);
\draw[ACliment]  (Cp) to node [right] {
} ($(C)+(0,.4)$);


\end{tikzpicture}
\caption{A second-order rewriting system.}
\label{FA7}
\end{center}
\end{figure}

There are also differences between our approach and that by polygraphs. It is worth noting that our objects of interest have a structure that is more complex than that of the objects considered by Burroni. In fact, the bottom row of Figure~\ref{FA7} consists of many-sorted $n$-categorical $\Sigma$-algebras.  Moreover, while Burroni focuses his attention in the presentation of free $n$-categories, our focus is in the presentation of many-sorted $n$-categorical $\Sigma$-algebras. Additionally, our sets of interest have an associated Artinian (pre)order that enables us to make proofs by induction and definitions by recursion on the different cells that make up the above $n$-categorical structure.

Another crucial difference lies in the fact that, in our approach, we work with translations and not with substitutions. Actually, the first part of our work, which deals with first-order structures and morphisms between them, could have been presented using substitutions instead of translations. But the second part, which deals with second-order structures, required a different treatment. In a second-order path it is essential that those parts of it that remain unchanged after applying a second-order production, kept their normal form to avoid duplications. Consequently, in order to unify the treatment, we had to use translations in the two parts. On the contrary, in the polygraph approach, we think that this choice is not necessary since a second-order path being a term representing formal compositions, those parts of it that remain unchanged will be, by default, in normal form.

In our approach, we also work with higher-order categories and face the same notational and conceptual difficulties as in the polygraph approach. However, as far as the representation of cells is concerned, our proposal advocates using either derivations or terms representing them interchangeably, leading to a Curry-Howard type result. We think that following the dictum that ``derivations are terms'' and its natural generalisation ``the more complex the derivations, the more complex the terms will be'', can be very useful. The representation of $n$-cells as terms has the virtue of making recursion between layers explicit  since, taking advantage of the representation of the $n$-cells as terms, an $(n+1)$-cell is a derivation between $n$-cells. That is, the process for layer $n$ is repeated for layer $n+1$, always by rewriting terms. This also applies to non-categorical operations, which adds to the difficulty of the study.
In this context, the many-sorted Curry-Howard mapping serves as a key tool to determine when two paths belong to the same class. This can be easily determined by checking if they have the same image under the many-sorted Curry-Howard mapping. In contrast, in the case of polygraphs, a quotient is taken relative to the smallest congruence containing the categorial axioms, which entails a loss of control over the classes and a shift towards diagrammatic models instead of term models. In addition, representing cells as terms allows for computational processing, which can be complicated with other more graphical representations.

Finally, our proposal includes the description of morphisms between higher-order rewriting systems, a matter already introduced in the case of polygraphs. To our knowledge, the first definition of polygraph morphism appears in Guiraud and Malbos~\cite[Subsection 1.4.2, p.~10]{GM09}.  For a comprehensive exploration of morphisms in higher-order polygraphs, interested readers are directed to the compendium~\cite{abgmmm23}. Since polygraphs only consider the higher-order category structure, morphisms are functors satisfying certain commutativity conditions between diagrams. In our case, in contrast to the case of the polygraphs, in addition to comparing categories, we also have to compare many-sorted signatures. For this purpose we consider derivors between signatures. The consequence of this is an increase in the complexity of the definition of the morphisms between higher-order rewriting systems. We hope that our work will contribute to a better understanding of these morphisms.

\section{A brief note on the origins of rewriting}

It is challenging to fully and fairly acknowledge all those who have contributed to the development of term rewriting, a field with deep roots in both mathematics and mathematical logic. To this end, and without claiming to be exhaustive, we recommend the following references~\cite{BKV03, rb89, CH06, Mul24, Pow89, Pow14, Sel06, ST00}. 
In what follows, we focus on further elaborating the contributions of Thue, Sch\"{o}nfinkel, Curry, and Church
while also referring, where appropriate, to subsequent developments by other authors relevant to the topic under consideration. 

Rewriting in mathematics can be traced back to the early 20th century, with the groundbreaking work of Thue between 1910 and 1914. In this regard, and to begin with, it is worth recalling B\"{u}chi's illuminating and insightful remark on p.~235 of~\cite{rb89}, which is particularly pertinent to the present discussion:
\begin{quotation}
In the past 25 years I have met several happy people telling me they had just discovered that terms are trees. I was happy too, when I discovered this. But soon \emph{I found it all in Thue}~\cite{Th10}. This has to be the first investigation, both of trees and of general combinatorial rules. \!\ldots Tree-manipulation rules came first, and these led him onto the much broader investigation of manipulation of strings of symbols (words) by the general rule $\xi a\eta\rightarrow\xi b\eta$ (see Thue~\cite{Th12, Th14}).
You should study these works, if you want to know how mathematical linguistics, equational logic, and the general idea of an algorithm were born. \emph{Where there seemed little room for mathematical investigation, Thue's clear mind created a science of strings of symbols and trees, showing the way to a rigorous theory of production rules and decidability}. In Thue~\cite{Th14} you will find what is now called ``the word problem for semigroups'' and in Thue~\cite{Th10} ``the word problem for free algebras''. Thue spoke of ``insurmontable difficulties'' that might exist. So he clearly understood that a decision problem may be unsolvable. I think he was the first to realize this [Tietze, in~\cite{T08}, said of finitely presented groups: ``The question whether two groups are isomorphic is not generally solvable'', we add]. \emph{(Emphasis ours; note that B\"{u}chi cites \cite{Th10, Th14}, but \cite{Th10} appears to be incorrect---we believe \cite{Th12, Th14} is intended.)}
\end{quotation}
B\"{u}chi's words about Thue seem to echo Noether's words about Dedekind: ``Everything is already in Dedekind''. Thue's work in~\cite{Th10}, \cite{Th12} and \cite{Th14} was far ahead of his time, with many of his results later independently rediscovered by researchers unaware of his original contributions. Unlike Dedekind, however, Thue did not have a Noether---someone who, contemporaneously, recognized the depth and significance of the contributions in these papers.


In line with the preceding discussion, although Thue is cited only once in~\cite{BKV03}, specifically in Exercise~3.1.3, on p.~62, and solely in connection with the \emph{Thue-Morse sequence}, it is noteworthy that Thue's works~\cite{Th10, Th12, Th14} on \emph{rewriting}---reprinted in~\cite{Th77}, published in 1977---are not mentioned in~\cite{BKV03}, published in 2003. This omission deserves attention, given the relevance of Thue's aforementioned contributions to the topic and the content of~\cite{BKV03}.

Next, drawing on the excellent paper by Steinby and Thomas~\cite{ST00}, which provides a detailed analysis of Thue's masterpiece~\cite{Th10}---whose title translates into English as \emph{The solution of a special case of a general logical problem}---as well as on B\"{u}chi's book~\cite{rb89}, and reiterating that---despite its significance---this particular contribution by Thue has been largely overlooked in the literature, we present a concise exposition of the essential aspects of this work.

At the beginning of~\cite{Th10}, and using modern terminology, Thue considers trees, i.e., essentially many-sorted terms in $\mathrm{T}_{\Sigma}((\varnothing)_{s\in S})$, the underlying $S$-sorted set of $\mathbf{T}_{\Sigma}((\varnothing)_{s\in S})$, the free $\Sigma$-algebra on $(\varnothing)_{s\in S}$, the $S$-sorted set which is constantly empty, for a set of sorts $S$ and an $S$-sorted signature $\Sigma$.
Following this, Thue formulates the problem of whether the equality of two trees can be derived from a given set of axioms. But, to simplify matters, Thue restricts himself to the single-sorted case, i.e., to single-sorted algebras, and, in addition, allows terms with variables. Now, Thue poses the word problem for free single-sorted algebras as: Can we decide whether an equation $P = Q$, with $P$, $Q\in \mathrm{T}_{\Sigma}(\{x_{i}\mid i\in n\})$, for some single-sorted signature $\Sigma$ and some finite set of variables $\{x_{i}\mid i\in n\}$, holds for all values of the variables assuming that all equations in a given finite subset $\mathcal{E}$ of $\mathrm{T}_{\Sigma}(\{x_{i}\mid i\in n\})^{2}$ hold for all values of the variables? Before Thue sheds new light on the word problem, he provides, according to~\cite{rb89}, the first rigorous (recursive) definition of a tree. More precisely, he recursively defines the set of unlabelled (binary) trees. Then he defines a labelled tree as a tree together with a mapping that assigns labels to the nodes of the unlabelled tree. The labels are of two kinds: variables and binary operation symbols. After that Thue recursively defines the identity of two unlabelled trees and then the identity of two labelled trees. Next, related to the word problem, Thue defines the tree replacement operation and for a finite subset $\mathcal{E}$ of $\mathrm{T}_{\Sigma}(\{x_{i}\mid i\in n\})^{2}$ defines, by using the replacement operation, the relation of congruence between two terms. Following this, Thue defines the relation of equality between two terms, which is what now is called the Thue congruence generated by $\mathcal{E}$, so that the word problem for a pair of terms is reduced to asking whether such a pair of terms belongs to the Thue congruence generated by some set $\mathcal{E}$. At this point, as remarked by B\"{u}chi~\cite{rb89}, Thue said that ``There may well be insurmountable obstacles to a general solution of the problem'' (``problem'' refers to \emph{the word problem for equational theories}, i.e., to decide whether a system of equations implies a given equation). Following this and for a single-sorted signature $\Sigma$ with only one binary operation symbol and a set of variables with only one variable $x$, Thue converts his system of equations $\mathcal{E}\subseteq \mathrm{T}_{\Sigma}(\{x\})^{2}$, which is now closed under substitutions and subject to satisfying three conditions related to the complexity of the terms of the equations, to an equivalent and complete (i.e., terminating and confluent) rewriting system $\mathcal{R}$ by orienting the equations; and $\mathcal{R}$ gives an algorithm for solving the word problem for $\mathcal{E}$. In the remainder of his work, Thue discusses proofs in normal forms, local confluence and the Church-Rosser property.

Before proceeding further and in connection with Thue's trees, we would like to highlight that in Chapter~3, Section~1 of~\cite{Her30}, Herbrand stated, for a many-sorted signature, the definition of the notion of an absolutely free many-sorted algebra on a many-sorted set of variables. 
Moreover, the semi-Thue rules, originally introduced by Thue in~\cite{Th12, Th14} as rewrite rules on words, form a special case of Post canonical rules (see~\cite{Post43,Post44}). Semi-Thue rules are, in turn, closely related to Turing machines (cf.~\cite{D82}).

To summarize the above, it is best to once again quote B\"{u}chi's words from~\cite{rb89}, on p.~220:
\begin{quotation}
It all started with trees, which became words, and all was closely knitted in with universal algebra and decidability. Where have the good days vanished to?
\end{quotation}

It is worth noting that Dehn formulated the word problem---along with the conjugacy and isomorphism problems---for groups in~\cite{MD11}, one year after Thue posed the word problem for free algebras in~\cite{Th10}; and that Thue, in turn, formulated the word problem for semigroups in~\cite{Th14}, three years after Dehn's formulation of the word problem for groups in~\cite{MD11}. In this regard, Magnus, one of Dehn's students, wrote in~\cite{CM82}, on p.~54:
\begin{quotation}
What appears to be incidental or, if one prefers, miraculous, is the fact that  independent of Dehn and independent of topology, a contemporary mathematician had begun to ask questions of the type of the word problem in combinatorial group theory,  but in an even \emph{more general and highly abstract setting} [the word problem for free algebras in~\cite{Th10} and the word problem for semigroups in~\cite{Th14}, we add]. We are referring to the work of Thue, who may be considered as the founder of a general theory of semigroups [and, on the basis of~\cite{Th10}, like the one who laid the foundations for many-sorted universal algebra, we add]. 
With one widely quoted exception, this work of his is largely forgotten nowadays. We do not know whether Dehn was influenced by Thue, and we have reasons to doubt it. 
We know that Dehn knew Thue personally, but only very superficially. Dehn mentioned Thue's work on occasion, observing that \emph{Thue's papers dealt with combinatorial problems} [Dehn was not a Noether for Thue, we add]. But he never used them, and indeed there is no known direct application of Thue's work to Dehn's group-theoretic problems. (\emph{Emphasis ours}.) 
\end{quotation}

Nearly half a century after the publication of Dehn's~\cite{MD11} paper on the word problem, Novikov~\cite{N55} proved that it is recursively unsolvable. For further discussion, see the book by Chandler and Magnus~\cite{CM82}, Miller\! III's  paper~\cite{Mi14} as well as the work of M\"{u}ller-Stach~\cite{Mul24}, all of which provide additional insights into the development of this topic. Moreover, Thue's word problem for semigroups was later studied by Post~\cite{Post47} (using Turing machines) and Markov~\cite{Mar47} (using Post's normal systems), who independently proved the undecidability of the word problem for semi-Thue systems. Their work established rewriting as a fundamental tool for investigating computability and decision problems. 
For studies on the concept of unsolvability, see~\cite{D82}, \cite{D04}, and especially the article by Post~\cite{Post04} within~\cite{D04}, in which he recounts---with notable intellectual honesty and at times with a touch of melancholy---his decade-long anticipation of the fundamental results of G\"{o}del~\cite{G31} (the first incompleteness theorem), Church~\cite{Chu36} (Church's thesis, $\lambda$-definability and existence of unsolvable problems) and Turing~\cite{Tur36} (the nature of computability, Turing machines and unsolvability). 

In mathematical logic, rewriting emerged in the late 1920s and early 1930s with Sch\"{o}nfinkel's function calculus, Curry's combinatory logic and Church's $\lambda$-calculus.

Combinatory logic, which is a most important example of a term rewriting system, started with Sch\"{o}nfinkel~\cite{Sch24}. The main aim of Sch\"{o}nfinkel was to eliminate the variables from the predicate calculus and this he achieved by developing a function calculus, functions being understood as rules of correspondence (or functions in intension in the terminology of Church~\cite{Chu41}), and showing that every formula of the predicate calculus can be expressed by means of the particular functions (combinators in the terminlogy of Curry~\cite{CurI30, CurII30}) $C$, $S$, $U$ and application. Sch\"{o}nfinkel's key idea was to extend the notion of function so that functions can be arguments and values of other functions. 
In~\cite{LS86}, Lambek and Scott summarised Sch\"{o}nfinkel's contributions as follows: He considers algebras $\mathbf{A} = (A, {}^{\wr}, I,K,S)$, where ${}^{\wr}$ (\emph{application}) is a binary operation, and $I$ (\emph{Identit\"{a}tsfunktion}), $K$ (\emph{Konstanzfunction}) and $S$ (\emph{Verschmelzungsfunktion}) are elements of $A$ satisfying the following equations: $(\mathbf{I})$ $I^{\wr}a = a$, $(\mathbf{K})$ $(K^{\wr}a)^{\wr}b = a$ and $(\mathbf{S})$ $((S^{\wr}f)^{\wr}g)^{\wr}c = (f^{\wr}c)^{\wr}(g^{\wr}c)$, for every $a$, $b$, $c$, $f$ and $g$ in $A$. These algebras are now known as Sch\"{o}nfinkel algebras, a term introduced by Lambek and Scott in~\cite{LS86}. Moreover, Sch\"{o}nfinkel established the functional completeness of these  algebras: For every polynomial $P$ in an indeterminate $x$ over a Sch\"{o}nfinkel algebra $\mathbf{A}$, i.e., for every $P\in A[x]$, where $A[x]$ is the underlying set of $\mathbf{A}[x]$, the coproduct in the category of Sch\"{o}nfinkel algebras of $\mathbf{A}$ and the free Sch\"{o}nfinkel algebra on $\{x\}$, there exists an $f\in A$ such that $P = f^{\wr}x$. The functional completeness  can be generalized to polynomials in $\mathbf{A}[x_{0},\ldots,x_{n-1}]$. Let us note that the application operation allows for self-application, i.e., for every $f\in A$, we have that $f^{\wr}f\in A$; and the functional completeness allows to define functions by self-reference.

In~\cite{CurI30, CurII30}, Curry, working partly independently, rediscovered results originally due to Sch\"{o}nfinkel (thus operating, in particular, without the use of  bound variables). He extended the defining equations $(\mathbf{I})$, $(\mathbf{K})$, $(\mathbf{S})$ of Sch\"{o}nfinkel algebras by adding five additional equations $(\mathrm{c}1)$, \ldots, $(\mathrm{c}5)$ (see~\cite{St72}). The resulting structures are now referred to as Curry algebras, a term introduced by Lambek and Scott in~\cite{LS86}. It was within the framework of these algebras that Curry reinforced the notion of functional completeness, in the following sense: For every polynomial $P$ in an indeterminate $x$ over a Curry algebra $\mathbf{A}$, there exists a \emph{unique} element $f\in A$ such that $P = f^{\wr}x$. The just stated functional completeness can be generalized to polynomials in $\mathbf{A}[x_{0},\ldots,x_{n-1}]$,  with the same uniqueness property holding. Moreover, for $\mathbf{T}_{\Sigma}(V)$, the (absolutely) free algebra of combinatory terms on a countably infinite set of variables $V$ and for the signature $\Sigma = \{I,K,S,{}^{\wr}\}$, Curry shows that the $\eta$-equality relation, denoted by $=_{\eta}$---which is the congruence on $\mathbf{T}_{\Sigma}(V)$ generated by the ordered pairs of combinatory terms appearing in $(\mathbf{I})$, $(\mathbf{K})$, $(\mathbf{S})$, $(\mathrm{c}1)$, \ldots, $(\mathrm{c}5)$---has the following property: (Principle of Extensionality) If $x$ is a variable not occurring in the combinatory terms $P$ and $Q$, and if $P^{\wr}x =_ {\eta} Q^{\wr}x$, then $P =_ {\eta} Q$. For a discussion of Curry's combinatory logic from the point of view of term rewriting systems---namely, treating the axioms as rewrite rules---, see~\cite{BKV03}.



In~\cite{ChuI32,ChuII32}, Church introduced the $\lambda$-calculus (conversion calculus), which, according to him, is a set of postulates for the foundation of logic and to avoid the paradoxes of mathematical logic and different of the method of Russell (type theory) or that of Zermelo (axiomatic set theory). However, in 1935 Kleene and Rosser proved, in particular, that Church's full system of logic in~\cite{ChuI32,ChuII32} was inconsistent. This resulted in Church's restriction to the pure part of his theory that was expounded in~\cite{Chu41}. The $\lambda$-calculus is a fundamental example of a higher-order term rewriting system (in the sense of~\cite{BKV03}), and, unlike Curry, Church employs bound variables, via the $\lambda$-abstraction operator $\lambda x$, where $x$ is a variable. However, there exists a close connection between the lambda calculus and combinatory logic. More precisely, the theory of combinators with $\eta$-equality is equivalent to the $\lambda$-calculus with $\beta\eta$-conversion. Furthermore, we emphasize that the translation of $\lambda$-terms into combinatory terms fundamentally relies on the functional completeness of combinatory logic. For a discussion of Church's $\lambda$-calculus see~\cite{BKV03}.

Let us recall that the partial recursive functions are definable in both Curry's combinatory logic and Church's $\lambda$-calculus, see~\cite{CH06}, and that this fact, in conjunction with the method of G\"{o}delization and the theorems of the $\lambda$-calculus itself, was used by Church in~\cite{Chu36} to prove, by diagonalization, that the problem of deciding for any two formulas of his $\lambda$-calculus whether they are convertible one into the other, is unsolvable. Moreover, for a detailed treatment of the relationships between typed combinatory logic, typed $\lambda$-calculus, and category theory, we refer the reader to~\cite{LS86}.

To conclude this subsection, it is also worth recalling that, in Chapter~1, Section~6 of \cite{Her30}, Herbrand provides an algebraic criterion for tautology checking, in which rewriting plays a central role. Before sketching how Herbrand proceeds in this regard, we recall---using modern terminology---that, in Chapter~1, Section~1 of \cite{Her30}, for the algebraic signature $\Lambda$ consisting of the operation symbols $\neg$ and $\vee$, and for a countably infinite set of variables 
$V = \{v_{n}\mid n\in \mathbb{N}\}$(implicit in Herbrand), he defines a \emph{proposition} as an element of the underlying set of $\mathbf{T}_{\Lambda}(\{v_{n}\mid n\in \mathbb{N}\})$, the free $\Lambda$-algebra on $\{v_{n}\mid n\in \mathbb{N}\}$. Then, in Chapter~1, Section~6 of \cite{Her30}, he begins by considering $\mathbf{F}_{2}$, the field of remainders modulo $2$, interpreting $0$ as truth and $1$ as falsehood. 
Next, for every $n\geq 1$ and every polynomial $P$ in $\mathbf{F}_{2}[v_{0},\ldots,v_{n-1}]$, Herbrand uses the notation $P \equiv 0\pmod{2}$ to express that the polynomial operation $P^{\mathbf{F}_{2}}$ from $\mathrm{F}_{2}^{n}$ to $\mathrm{F}_{2}$, canonically associated to $P$, is identically zero, i.e., it evaluates to zero for all inputs in $\mathrm{F}_{2}^n$. (We recall that there exists a canonical surjective homomorphism $(\cdot)^{\mathbf{F}_{2}}$ from $\mathbf{F}_{2}[v_{0},\ldots,v_{n-1}]$ to $\mathbf{F}_{2}^{\mathrm{F}_{2}^{n}}$ that sends a polynomial $P$ to the polynomial operation $P^{\mathbf{F}_{2}}$, thus $P \equiv 0\pmod{2}$ means that $P\in \mathrm{Ker}((\cdot)^{\mathbf{F}_{2}})$.) 
Following this, for every $n\geq 1$, he assigns to each proposition 
$\varphi$ in the underlying set of $\mathbf{T}_{\Lambda}(\{v_{i}\mid i\in n\})$, the 
free $\Lambda$-algebra on $\{v_{i}\mid i\in n\}$, a polynomial $\overline{\varphi}$ in $\mathbf{F}_{2}[v_{0},\ldots,v_{n-1}]$ as follows: If $\varphi$ is $v_{i}$, for some $i\in n$, then $\overline{v_{i}} = v_{i}$; if $\varphi$ is $\neg \psi$, for some proposition $\psi$, then $\overline{\neg \psi} = 1+\overline{\psi}$; and if $\varphi$ is $\psi\vee \xi$, for some propositions $\psi$ and $\xi$, then 
$\overline{\psi\vee\xi} = \overline{\psi}\cdot\overline{\xi}$. 
(Thus, Herbrand defines a homomorphism $\overline{(\cdot)}$ from $\mathbf{T}_{\Lambda}(\{v_{i}\mid i\in n\})$ to a derived $\Lambda$-algebra of $\mathbf{F}_{2}[v_{0},\ldots,v_{n-1}]$.) 
After this, Herbrand proves that a proposition $\varphi$ is true (a tautology) if, and only if, $\overline{\varphi} \equiv 0\pmod{2}$, i.e., if, and only if, $\overline{\varphi}\in \mathrm{Ker}((\cdot)^{\mathbf{F}_{2}})$. (Herbrand's proof is by algebraic induction.) 
Then Herbrand, in response to the question: ``How can it be determined  whether a polynomial is always $0\pmod{2}$?'', answers as follows: For every polynomial $P$ in $\mathbf{F}_{2}[v_{0},\ldots,v_{n-1}]$, $P \equiv 0\pmod{2}$ if, and only if, all of its coefficients are even.  
(Herbrand's proof proceeds by induction on the number of variables in the polynomials, expressing a polynomial $P$ in $n+1$ variables as $P = \binom{v_{n}}{0}^{\sharp}(P) + v_{n}(\binom{v_{n}}{0}^{\sharp}(P)+ \binom{v_{n}}{1}^{\sharp}(P))$, where, for $k\in \mathrm{F}_{2}$, $\binom{v_{n}}{k}^{\sharp}$ is the canonical extension---up to $\mathbf{F}_{2}[v_{0},\ldots,v_{n}]$---of the mapping 
$\binom{v_{n}}{k}$ from $\{v_{i}\mid i\in n + 1\}$ to $\mathrm{F}_{2}[v_{0},\ldots,v_{n-1}]$ that sends each $v_{i}$, for $i\in n$, to $v_{i}$ itself, and $v_{n}$ to $k$.)    
Finally, Herbrand states that, to verify whether a proposition $\varphi$ is true, one replaces it by $\overline{\varphi}$, and seeks to reduce $\overline{\varphi}$---i.e., \emph{rewrite} it to $0$---by using the (\emph{rewrite}) rules $x^{n}\equiv x \pmod{2}$ (for $n\geq 1$) and $1+1\equiv 0 \pmod{2}$.


\section{Prerequisites, notation and terminology}

The only prerequisites for reading this work are familiarity with category theory, universal algebra, in particular many-sorted partial algebras, the theory of ordered sets and set theory. This can be achieved by reading, for category theory, e.g., \cite{hs73,sM98}, for universal algebra, e.g., \cite{gb15,pb82,pb86,pb02,bs81,gtw85,gG08,agk67,sch66,sch70,w92}, for the theory of ordered sets, e.g., \cite{bl05,DP02,R08}, and, for set theory, e.g., \cite{nB70,end77}.

Our underlying set theory is $\mathbf{ZFSk}$, Zermelo-Fraenkel-Skolem set theory (also known as $\mathbf{ZFC}$, i.e.,  Zermelo-Fraenkel set theory with the axiom of choice) plus the existence of a Grothendieck universe $\boldsymbol{\mathcal{U}}$, fixed once and for all (see~\cite{sM98}, pp.~21--24). We recall that the elements of $\boldsymbol{\mathcal{U}}$ are called $\boldsymbol{\mathcal{U}}$-small sets and the subsets of $\boldsymbol{\mathcal{U}}$ are called $\boldsymbol{\mathcal{U}}$-large sets or classes. Moreover, from now on $\mathsf{Set}$ stands for the category of sets, i.e., the category whose set of objects is $\boldsymbol{\mathcal{U}}$ and whose set of morphisms is the set of all mappings between $\boldsymbol{\mathcal{U}}$-small sets.

In all that follows we use standard concepts and constructions from category theory, universal algebra, the theory of ordered sets and set theory. 
Nevertheless,
regarding set theory, we have adopted the following conventions.
An \emph{ordinal} $\alpha$ is a transitive set that is well-ordered by $\in$, thus $\alpha = \{\,\beta\mid \beta\in \alpha\,\}$. The first transfinite ordinal $\omega_{0}$ will be denoted by  $\mathbb{N}$, which is the set of all \emph{natural numbers}, and, from what we have just said about the ordinals, for every $n\in \mathbb{N}$, $n = \{0, \ldots,n-1\}$. If $\Phi$ and $\Psi$ are (binary) relations in a set $A$, then we will say that $\Psi$ is a \emph{refinement} of $\Phi$ if $\Psi\subseteq \Phi$.
We will denote by $\mathrm{Pfnc}(A,B)$ the set of all \emph{partial functions} from $A$ to $B$, and by $\mathrm{Fnc}(A,B)$ the set of all \emph{functions} from $A$ to $B$. We recall that a \emph{partial function} from $A$ to $B$ is a subset $F$ of $A\times B$ such that, for every $x\in A$, there is at most one $y\in B$ such that $(x,y)\in F$, and that a \emph{function} from $A$ to $B$ is a subset $F$ of $A\times B$ such that, for every $x\in A$, there exists a unique $y\in B$ such that $(x,y)\in F$. A partial function from $A$ to $B$ is usually denoted by $(F_{x})_{x\in \mathrm{Dom}(F)}$, where $\mathrm{Dom}(F)$, the \emph{domain (of definition)} of $F$, is $\{x\in A\mid \exists\,y\in B\,((x,y)\in F)\}$, and a function $F$ from $A$ to $B$ is generally denoted by $(F_{x})_{x\in A}$. We will denote by $\mathrm{Hom}_{\mathrm{p}}(A,B)$ the set of all \emph{partial mappings} from $A$ to $B$, and by $\mathrm{Hom}(A,B)$ (and, sometimes, also by $B^{A}$) the set of all \emph{mappings} from $A$ to $B$. We recall that a \emph{partial mapping} from $A$ to $B$ is an ordered triple $f = (A,F,B)$, denoted by $f\colon A\dmor B$, in which $F$ is a partial function from $A$ to $B$, and that a \emph{mapping} from $A$ to $B$ is an ordered triple $f = (A,F,B)$, denoted by $f\colon A\mor B$, in which $F$ is a function from $A$ to $B$. Therefore $\mathrm{Hom}_{\mathrm{p}}(A,B) = \{A\}\times \mathrm{Pfnc}(A,B)\times\{B\}$ and $\mathrm{Hom}(A,B) = \{A\}\times \mathrm{Fnc}(A,B)\times\{B\}$. For a partial mapping $f\colon A\dmor B$ we will denote by $\Gamma_{\!f}$ the underlying partial function of $f$ and by $\mathrm{Dom}(f)$ precisely $\mathrm{Dom}(\Gamma_{\!f})$. 
Let us note that, for a partial mapping $f\colon A\dmor B$, $\mathrm{d}_{0}(f)$, the categorial domain (or source) of $f$, which is $A$, contains $\mathrm{Dom}(f)$, while, for a mapping $f\colon A\mor B$, $\mathrm{d}_{0}(f) = \mathrm{Dom}(f)$. We will denote by $\mathrm{Sub}(A)$ the set of all sets $X$ such that $X\subseteq A$ and if $X\in\mathrm{Sub}(A)$, then we will denote by $\complement_{A}X$ or $A-X$ the complement of $X$ in $A$. Moreover, if $f$ is a mapping from $A$ to $B$, then the mapping $f[\cdot]$ from $\mathrm{Sub}(A)$ to $\mathrm{Sub}(B)$, of $f$-\emph{direct image formation}, sends $X$ in $\mathrm{Sub}(A)$ to
$
f[X] = \{y\in B\mid \exists\,x\in X\,(y = f(x))\}
$
in $\mathrm{Sub}(B)$, and the mapping $f^{-1}[\cdot]$ from $\mathrm{Sub}(B)$ to $\mathrm{Sub}(A)$, of $f$-\emph{inverse image formation}, sends $Y$ in $\mathrm{Sub}(B)$ to
$
f^{-1}[Y] = \{x\in A\mid f(x)\in Y\}
$
in $\mathrm{Sub}(A)$. In the sequel, for a mapping $f$ from $A$ to $B$ and a subset $X$ of $A$, we will write $\mathrm{Ker}(f)$ for the kernel of $f$, $\mathrm{Im}(f)$ for the image of $f$, to mean $f[A]$, and the restriction of $f$ to $X$ will be denoted by $f\!\!\upharpoonright_{X}$. Let us note that for partial mappings $f\colon A\dmor B$ and $g\colon B\dmor C$, the domain of the composite $g\circ f$ is $\mathrm{Dom}(g\circ f) = f^{-1}[\mathrm{Dom}(g)\cap f[\mathrm{Dom}(f)]]$ and the image is $\mathrm{Im}(g\circ f) = g[\mathrm{Dom}(g)\cap f[\mathrm{Dom}(f)]]$. 

More specific assumptions, conditions, and conventions will be included in subsequent chapters.

\section{A note from the authors}
This project started in the fall of 2018 under the guidance of Juan Climent and Enric Cosme.  Ra\'{u}l Ruiz joined our team in October 2023, and his contributions  constitute an important part of his PhD thesis.

On February 2024, we publicly unveiled our first findings on ArXiv.  This delay in the public presentation of our results was due, on the one hand, to the fact that we wanted to make the presentation as clean as possible and, on the other hand, to the fact that we wanted to present the results in such a way that the relationship between the results of the first and the second part was seen as inevitable and fundamental.  On January 2026, we published the second version on ArXiv, which includes the third part of this work, on morphisms of rewriting systems.

Due to the meticulous and conceptually intricate nature of our project, which has necessitated numerous revisions to achieve a state satisfactory to its authors, the document has expanded beyond our initial expectations. We hope that the latter will be offset by the former.

\section{Acknowledgements}
The second and third authors were supported by the  PID2024-159495NB-I00 grant from the
Ministerio de Ciencia e Innovaci\'{o}n, Spain, and the CIAICO/2023/007 grant  from the
Conselleria d’Educaci\'{o}, Universitats i Ocupaci\'{o}, Generalitat Valenciana. The second author has held a Specially Appointed Professor position at Nantong University during the completion of this work. The third author was supported by the  CIACIF/2022/489 grant from the Conselleria d'Educaci\'{o}, Universitats i Ocupaci\'{o}, Generalitat Valenciana co-funded by the European Social Fund. Special thanks to Asunci\'{o}n de Montesa for her categorical knowledge. 

\chapter{Many-sorted sets and many-sorted algebras}\label{S0B}

In this section we collect the basic facts, without proofs, about many-sorted sets and many-sorted (total) algebras, with special emphasis on the free many-sorted algebra  that we will need afterwards.

\section{Many-sorted sets}

\begin{assumption}
From now on $S$ stands for a set of sorts in $\boldsymbol{\mathcal{U}}$, fixed once and for all.
\end{assumption}

\begin{definition}
An $S$-\emph{sorted set}\index{many-sorted!set} is a mapping $A = (A_{s})_{s\in S}$ from $S$ to $\boldsymbol{\mathcal{U}}$. If $A$ and $X$ are $S$-sorted sets, then we will say that $X$ is an 
$S$-\emph{sorted subset} or, to abbreviate, a \emph{subset}\index{many-sorted!subset} of $A$, denoted by $X\subseteq A$, if, for every $s\in S$, $X_{s}\subseteq A_{s}$. We will denote by $\mathrm{Sub}(A)$ the set of all $S$-sorted sets $X$ such that $X\subseteq A$.
\end{definition}

\begin{definition}
Let $A$ and $B$ be $S$-sorted sets.  The \emph{cartesian product of} $A$ \emph{and} $B$, denoted by $A\times B$, is the $S$-sorted set $\left(A_{s}\times B_{s}\right)_{s\in S}.$

Let $\Phi$ be an $S$-sorted set. We will say that $\Phi$ is an $S$-\emph{sorted relation from}\index{many-sorted!relation} $A$ \emph{to} $B$ if $\Phi\subseteq A\times B$. Thus, for every $s\in S$, $\Phi_{s}\subseteq A_{s}\times B_{s}$. We denote by $\Rel(A,B)$ the set of all $S$-sorted relations from $A$ to $B$. If $A=B$, then we write $\Rel(A)$ instead of $\Rel(A,A)$ and call its elements $S$-\emph{sorted relations on} $A$. 

The \emph{diagonal of} $A$, denoted by $\Delta_{A}$, is the $S$-sorted relation on $A$ defined, for every $s\in S$, as $\Delta_{A_{s}}$, and the \emph{codiagonal of} $A$, denoted by $\nabla_{A}$, is the $S$-sorted relation on $A$ defined, for every $s\in S$, as $\nabla_{A_{s}} = A_{s}\times A_{s}$.

Let $\Phi$ be an $S$-sorted relation from $A$ to $B$ and $\Psi$ an $S$-sorted relation from $B$ to $C$, then the \emph{composition of} $\Phi$ \emph{and} $\Psi$, denoted by $\Psi\comp\Phi$, is the $S$-sorted relation from $A$ to $C$ defined, for every $s\in S$, as 
$$
\Psi_{s}\comp\Phi_{s}
=\left\lbrace
(x,z)\in A_{s}\times C_{s}
\mid \exists\,y\in B_{s}\,
\left(
(x,y)\in \Phi_{s}\And (y,z)\in \Psi_{s}
\right) 
\right\rbrace.
$$ 
This composition is associative and the diagonal relation $\Delta_{A}$ is a neutral element for it. 

Let $\Phi$ an $S$-sorted relation from $A$ to $B$. Then the \emph{inverse} of $\Phi$, denoted by $\Phi^{-1}$, is the $S$-sorted relation from $B$ to $A$ defined, for every $s\in S$, as
$$
\Phi_{s}^{-1}
=\left\lbrace(y,x)\in B_{s}\times A_{s}\mid (x,y)\in \Phi_{s}
\right\rbrace.
$$

Let $\Phi$ be an $S$-sorted relation from $A$ to $B$, $X$ an $S$-sorted subset of $A$, and $Y$ an $S$-sorted
subset of $B$, then the \emph{direct image} 
of $X$ under $\Phi$, denoted by $\Phi[X]$, is the $S$-sorted subset of $B$ defined, for every $s\in S$, as
$$
\Phi[X]_{s} =\left\lbrace\,b\in B_ {s}\mid \exists\, x\in X_ {s}((x,b)\in \Phi_
{s})\,\right\rbrace,
$$
therefore $\Phi[X] = (\Phi_ {s}[X_ {s}])_ {s\in S}$, and the \emph{inverse image} of $Y$ under $\Phi$, denoted by $\Phi^{-1}[Y]$, is the $S$-sorted subset of $A$ defined, for every $s\in S$, as
$$
\Phi^{-1}[Y]_{s} =\left\lbrace\,a\in A_ {s}\mid \exists\, y\in Y_ {s}((a,y)\in \Phi_
{s})\,\right\rbrace,
$$
therefore $\Phi^{-1}[Y] = (\Phi^{-1}_{s}[Y_ {s}])_{s\in S}$. For an $S$-sorted relation $\Phi$ on an $S$-sorted set $A$ we will call $\Phi[A]$ and $\Phi^{-1}[A]$ the \emph{image} and the \emph{domain} of $\Phi$ respectively and we will denote them by $\mathrm{Im}(\Phi)$ and $\mathrm{Dom}(\Phi)$ respectively. Moreover, we will call 
$\mathrm{Dom}(\Phi)\cup\mathrm{Im}(\Phi)$ the \emph{field} of $\Phi$ and we will denote it by $\mathrm{Fld}(\Phi)$.
    
An $S$-\emph{sorted function from}\index{many-sorted!function} $A$ \emph{to} $B$ is a functional $S$-sorted relation $F$ from $A$ to $B$, i.e., an $S$-sorted relation $F$ from $A$ to $B$ such that, for every $s\in S$, $F_{s}$ is a function from $A_{s}$ to $B_{s}$. We denote by $\Fnc(A,B)$ the set of all $S$-sorted functions from $A$ to $B$.  The composition of $S$-sorted functions, which is a particular case of the composition of $S$-sorted relations, is an $S$-sorted function.
      
An $S$-\emph{sorted mapping from}\index{many-sorted!mapping} $A$ \emph{to} $B$ is a triple $f = (A,F,B)$ where
$F$ is an $S$-sorted function from $A$ to $B$. We denote by $\Hom(A,B)$ or by $B_{A}$ the set of all
$S$-sorted mappings from $A$ to $B$. We consider the expressions $f\in\Hom(A,B)$, $f\in B_{A}$, and $f\colon A\mor B$ as synonymous. Moreover, given $f\colon A\mor B$ and $g\colon B\mor C$, $g\comp f = (A,G\comp F,C)$,
the \emph{composition of} $f$ \emph{and} $g$, is an $S$-sorted mapping from $A$ to $C$, and $\mathrm{id}^{A} = (A,\Delta_{A},A)$, is an  $S$-sorted endomapping of $A$, the identity $S$-sorted mapping at $A$. We denote by
$\mathsf{Set}^{S}$ the category of $S$-sorted sets and $S$-sorted mappings.

Let $f\colon A\mor B$ be an $S$-sorted mapping. Then the mapping 
$$
f[\cdot]\colon\mathrm{Sub}(A)\mor\mathrm{Sub}(B),
$$ 
of $f$-\emph{direct image formation}, sends $X \in \mathrm{Sub}(A)$ to $f[X] = (f_{s}[X_{s}])_{s\in S} \in \mathrm{Sub}(B)$, and the mapping 
$$
f^{-1}[\cdot]\colon\mathrm{Sub}(B)\mor\mathrm{Sub}(A),
$$ 
of $f$-\emph{inverse image formation}, sends $Y \in \mathrm{Sub}(B)$ to $f^{-1}[Y] = (f_{s}^{-1}[Y_{s}] )_{s\in S} \in \mathrm{Sub}(A)$. 

The \emph{image} of $f$, denoted by $\mathrm{Im}(f)$, is $f[A]$, i.e., $(f_{s}[A_{s}])_{s\in S}$. Moreover, if $X\subseteq A$, then the \emph{restriction of} $f$ \emph{to} $X$, denoted by $f\!\!\upharpoonright_{X}$, $f|_{X}$ or $\mathrm{res}_{X}(f)$, is $f\circ \mathrm{in}^{X,A}$, where $\mathrm{in}^{X,A} = (\mathrm{in}^{X_{s},A_{s}})_{s\in S}$ is the canonical embedding of $X$ into $A$.
\end{definition}

\begin{remark}
To give an $S$-sorted mapping from $A$ to $B$, as in the just stated definition, is equivalent to give an $S$-indexed family $f = (f_{s})_{s\in S}$, where, for every $s$ in $S$, $f_{s}$ is a mapping from $A_{s}$ to  $B_{s}$. Thus, an $S$-sorted mapping from $A$ to $B$ is, essentially, an element of $\prod_{s\in S}\mathrm{Hom}(A_{s}, B_{s})$. 
\end{remark}

\begin{proposition}\label{PDeltaPhiFunc}
Let $S$ and $T$ be sets and $\varphi \colon S \mor T$ a mapping. We let $\Delta_{\varphi}$ stand for the assignment from $\mathsf{Set}^{T}$ to $\mathsf{Set}^{S}$ defined as follows:
\begin{enumerate}
\item
for every $T$-sorted set $A = (A_{t})_{t \in T}$, $\Delta_{\varphi}(A)$ is the $S$-sorted set $A_{\varphi} = (A_{\varphi(s)})_{s \in S}$, and
\item
for every $T$-sorted sets $A$ and $B$, and every $T$-sorted mapping $f \colon A \mor B$, $\Delta_{\varphi}(f)$ is the $S$-sorted mapping $f_{\varphi} = (f_{\varphi(s)})_{s \in S}$.
\end{enumerate}
Then, $\Delta_{\varphi}$ is a covariant functor from $\mathsf{Set}^{T}$ to $\mathsf{Set}^{S}$.
\end{proposition}

\begin{proposition}\label{PMSetFunc}
We let $\mathrm{MSet}$ stand for the assignment from $\mathsf{Set}$ to $\mathsf{Cat}$ defined as follows:
\begin{enumerate}
\item
for every set $S$, $\mathrm{MSet}(S)$ is the category $\mathsf{Set}^{S}$, and
\item
for every mapping $\varphi \colon S \mor T$, $\mathrm{MSet}(f)$ is $\Delta_{\varphi}$ the functor from $\mathsf{Set}^{T}$ to $\mathsf{Set}^{S}$ introduced in Proposition~\ref{PDeltaPhiFunc}.
\end{enumerate}
Then $\mathrm{MSet}$ is a contravariant functor from $\mathsf{Set}$ to $\mathsf{Cat}$.
\end{proposition}

\begin{definition}
Let $I$ be a set in $\boldsymbol{\mathcal{U}}$ and $(A^{i})_{i\in I}$ an $I$-indexed family of $S$-sorted sets. Then the \emph{product}\index{many-sorted!product} of $(A^{i})_{i\in I}$, denoted by $\prod_{i\in I}A^{i}$, is the $S$-sorted set defined, for every $s\in S$, as $\left(\prod\nolimits_{i\in I}A^{i}\right)_{s} = \prod\nolimits_{i\in I}A^{i}_{s}$, where
$$
\textstyle
\prod_{i\in I}A^{i}_{s} = \left\{(a_{i})_{i\in I}\in\mathrm{Fnc}\left(I,\bigcup_{i\in I}A^{i}_{s}\right)\mid \forall\,i\in I\,\left(a_{i}\in A^{i}_{s}\right)\right\}.
$$
For every $i\in I$, the \emph{i-th canonical projection}, $\mathrm{pr}^{i} = (\mathrm{pr}^{i}_{s})_{s\in S}$, is the $S$-sorted mapping from  $\prod_{i\in I}A^{i}$ to $A^{i}$ that, for every $s\in S$, sends $(a_{i})_{i\in I}$ in $\prod_{i\in I}A^{i}_{s}$ to $a_{i}$ in $A^{i}_{s}$. The ordered pair $(\prod_{i\in I}A^{i},(\mathrm{pr}^{i})_{i\in I})$ has the following universal property: For every $S$-sorted set $B$ and every $I$-indexed family of $S$-sorted mappings $(f^{i})_{i\in I}$, where, for every $i\in I$, $f^{i}$ is an $S$-sorted mapping from $B$ to $A^{i}$, there exists a unique $S$-sorted mapping $\left<f^{i}\right>_{i\in I}$ from $B$ to $\prod_{i\in I}A^{i}$ such that, for every $i\in I$, $\mathrm{pr}^{i}\circ \left<f^{i}\right>_{i\in I} = f^{i}$.

The \emph{coproduct}\index{many-sorted!coproduct} of $(A^{i})_{i\in I}$, denoted by $\coprod_{i\in I}A^{i}$, is the $S$-sorted set defined, for every $s\in S$, as $\left(\coprod\nolimits_{i\in I}A^{i}\right)_{s} = \coprod\nolimits_{i\in I}A^{i}_{s}$, where
$$
\textstyle
\coprod_{i\in I}A^{i}_{s} = \bigcup_{i\in I}(A^{i}_{s}\times\{i\}).
$$
For every $i\in I$, the \emph{i-th canonical injection}, $\mathrm{in}^{i}$, is the $S$-sorted mapping from $A^{i}$ to $\coprod_{i\in I}A^{i}$ that, for every $s\in S$, sends $a$ in $A^{i}_{s}$ to $(a,i)$ in $\coprod_{i\in I}A^{i}_{s}$. The ordered pair $(\coprod_{i\in I}A^{i},(\mathrm{in}^{i})_{i\in I})$ has the following universal property: For every $S$-sorted  set $B$ and every $I$-indexed family of $S$-sorted mappings $(f^{i})_{i\in I}$, where, for every $i\in I$, $f^{i}$ is an $S$-sorted mapping from $A^{i}$ to $B$, there exists a unique $S$-sorted mapping $[f^{i}]_{i\in I}$ from $\coprod_{i\in I}A^{i}$ to $B$ such that, for every $i\in I$, $ [f^{i}]_{i\in I}\circ\mathrm{in}^{i} = f^{i}$.

The remaining set-theoretic operations on $S$-sorted sets: $\amalg$ (binary coproduct), $\bigcup$ (union), $\cup$ (binary union), $\bigcap$ (intersection), $\cap$ (binary intersection), $-$ (difference), and $\complement_{A}$ (complement of an $S$-sorted set in a fixed $S$-sorted $A$), are defined in a similar way, i.e., componentwise.
\end{definition}

\begin{definition}
We will denote by $1^{S}$ the (standard) final $S$-sorted set of $\mathsf{Set}^{S}$, which is $1^{S} = (1)_{s\in S}$, and by $\varnothing^{S}$ the initial $S$-sorted set, which is $\varnothing^{S} = (\varnothing)_{s\in S}$. We shall abbreviate $1^{S}$ to $1$ and $\varnothing^{S}$ to $\varnothing$ when this is unlikely to cause confusion.
\end{definition}


\begin{definition}
Let $\delta$ be the mapping from $S\times \boldsymbol{\mathcal{U}}$ to $\boldsymbol{\mathcal{U}}^{S}$ that sends $(t,X)$ in $S\times \boldsymbol{\mathcal{U}}$ to the $S$-sorted set $\delta^{t,X} = (\delta^{t,X}_{s})_{s\in S}$ defined, for every $s\in S$, as follows: $\delta^{t,X}_{s} = X$, if $s = t$; $\delta^{t,X}_{s} = \varnothing$, otherwise. We will call the value of $\delta$ at $(t,X)$ the \emph{delta of Kronecker associated to}\index{Delta of Kronecker} $(t,X)$. If $X = \{x\}$, then, for simplicity of notation, we will write $\delta^{t,x}$ instead of $\delta^{t,\{x\}}$. Moreover, for a sort $t$ in $S$, $\delta^{t,1}$, the delta of Kronecker associated to $(t,1)$, will be denoted by $\delta^{t}$ and called \emph{delta of Kronecker}.

%
\end{definition}

\begin{remark}
For a sort $t\in S$ and a set $X$, the $S$-sorted set $\delta^{t,X}$ is isomorphic to the $S$-sorted set $\coprod_{x\in X}\delta^{t}$, i.e., to the coproduct of the family $(\delta^{t})_{x\in X}$ which is constantly $\delta^{t}$.

For every sort $t\in S$ we have a functor $\delta^{t,\cdot}$ from $\mathsf{Set}$ to $\mathsf{Set}^{S}$. In fact, for every set $X$, $\delta^{t,\cdot}(X) = \delta^{t,X}$, and, for every mapping $f\colon X\mor Y$, $\delta^{t,\cdot}(f) = \delta^{t,f}$, where, for $s\in S$, $\delta^{t,f}_{s} = \mathrm{id}^{\varnothing}$, if $s\neq t$, and $\delta^{t,f}_{t} = f$. Moreover, for every $t\in S$, the object mapping of the functor $\delta^{t,\cdot}$ is injective and $\delta^{t,\cdot}$ is full and faithful. Hence, for every $t\in S$, $\delta^{t,\cdot}$ is a full embedding from $\mathsf{Set}$ to $\mathsf{Set}^{S}$.

The final object $1^{S}$ does not generate ($\equiv$ separate) the category $\mathsf{Set}^{S}$. However, the set $\{\,\delta^{s}\mid s\in S\,\}$, of the deltas of Kronecker, is a generating ($\equiv$ separating) set for the category $\mathsf{Set}^{S}$. Therefore, every $S$-sorted set $A$ can be represented as a coproduct of copowers of deltas of Kronecker, i.e., $A$ is naturally isomorphic to $\coprod_{s\in S}\mathrm{card}(A_{s})\boldsymbol{\cdot}\delta^{s}$, where, for every $s\in S$, $\mathrm{card}(A_{s})\boldsymbol{\cdot}\delta^{s}$ is the copower of the family $(\delta^{s})_{\alpha\in \mathrm{card}(A_{s})}$, i.e., the coproduct of $(\delta^{s})_{\alpha\in \mathrm{card}(A_{s})}$.
To this we add the following facts: (1) $\{\,\delta^{s}\mid s\in S\,\}$ is the set of atoms of the Boolean algebra $\mathbf{Sub}(1^{S})$, of subobjects of $1^{S}$; (2) the Boolean algebras $\mathbf{Sub}(1^{S})$ and $\mathbf{Sub}(S)$ are isomorphic; (3) for every $s\in S$, $\delta^{s}$ is a projective object; and (4) for every $s\in S$, every $S$-sorted mapping from $\delta^{s}$ to another $S$-sorted set is a monomorphism.

In view of the above, it must  be concluded that the deltas of Kronecker are of crucial importance for many-sorted sets and associated fields.
\end{remark}

Before proceeding any further, let us point out that it is no longer unusual to find in the works devoted to investigate both many-sorted algebras and many-sorted algebraic systems the following. (1) That an $S$-sorted set $A$ is defined in such a way that $\mathrm{Hom}(1^{S},A)\neq \varnothing$, or, what is equivalent, requiring that, for every $s\in S$, $A_{s}\neq \varnothing$. This has as an immediate consequence that the corresponding category is not even finite cocomplete. Since cocompleteness (and completeness) are desirable properties for a category, we exclude such a convention in our work (the admission of $\varnothing^{S}$ is crucial in many applications). And (2) that an $S$-sorted set $A$ must be such that, for every $s$, $t\in S$, if $s\neq t$, then $A_{s}\cap A_{t} = \varnothing$. We also exclude such a requirement (the possibility of a common underlying set for the different sorts is very important in many applications). The above conventions are possibly based on the untrue widespread belief that many-sorted equational logic and many-sorted first-order logic with equality are \emph{inessential} variations of equational logic and first-order logic with equality, respectively. One can find a definitive refutation to the just mentioned belief in~\cite{gm85} and \cite{m76}, regarding many-sorted equational logic, and in~\cite{Hook85}, with respect to many-sorted first-order logic with equality.

\begin{figure}
$$\xymatrix{
  & A \ar[dl]_{h} \ar[d]^{f} \\
Y \ar @{{ +}{-}{>}}[r]_{\mathrm{in}^{Y,B}} & B
 }
$$
\caption{The corestriction of a many-sorted mapping.}
\label{FPUS}
\end{figure}

\begin{proposition}\label{PUS}
Let $B$ be an $S$-sorted set, $Y$ a subset of $B$, and $f$ an $S$-sorted mapping from $A$ to $B$. Then the following statements are equivalent:
\begin{enumerate}
\item $\mathrm{Im}(f)\subseteq Y$.
\item There exists an $S$-sorted mapping $h$ from $A$ to $Y$ such that the diagram
in Figure~\ref{FPUS} 
commutes.
\end{enumerate}
If one of the above equivalent statements holds, then we will call $h$, which is univocally determined, the \emph{corestriction of} $f$ \emph{to} $Y$ and we denote it by $f|^{Y}$ or $\mathrm{cores}_{Y}(f)$.
\end{proposition}


\begin{proposition}
Let $\mathsf{PSet}^{S}$ be the category whose objects are the \emph{$S$-sorted set pairs}, i.e., the ordered pairs pairs $(A,X)$ where $A$ is an $S$-sorted set and $X\subseteq A$, and in which the set of morphisms from $(A,X)$ to $(B,Y)$ is the set of all $S$-sorted mappings $f$ from $A$ to $B$ such that $f[X]\subseteq Y$. Let $G$ be the functor from $\mathsf{Set}^{S}$ to $\mathsf{PSet}^{S}$ whose object mapping sends $A$ to 
$(A,A)$ and whose morphism mapping sends $f\colon A\mor B$ to $f\colon (A,A)\mor (B,B)$. Then, for every $S$-sorted pair $(B,Y)$, there exists a universal mapping from $G$ to $(B,Y)$, which is precisely the ordered pair $(Y,\mathrm{in}^{Y,B})$ with $\mathrm{in}^{Y,B}\colon (Y,Y)\mor (B,Y)$ the morphism in $\mathsf{PSet}^{S}$ associated to $\mathrm{in}^{Y,B}\colon Y\mor B$.
\end{proposition}

\begin{definition}
Let $f$ be a $\mathsf{PSet}^{S}$-morphism from $(A,X)$ to $(B,Y)$. Then we denote by $f|_{X}^{Y}$, $\mathrm{bires}_{X,Y}(f)$, or, if no confusion can arise, $\widehat{f}$ the $S$-sorted mapping $\mathrm{cores}_{Y}(\mathrm{res}_{X}(f))$ (which is identical to $\mathrm{res}_{X}(\mathrm{cores}_{Y}(f))$). We will call this $S$-sorted mapping the \emph{birestriction of} $f$ to $X$ and $Y$. 
\end{definition}

\begin{definition}
Let $A$ be an $S$-sorted set. Then the \emph{cardinal}\index{many-sorted!cardinal} of $A$, denoted by $\mathrm{card}(A)$, is $\mathrm{card}(\coprod A)$, i.e., the cardinal of the set $\coprod A = \bigcup_{s\in S}(A_{s}\times \{s\})$. An $S$-sorted set $A$ is \emph{finite}\index{many-sorted!finite} if $\mathrm{card}(A)<\aleph_{0}$. We will say that an $S$-sorted set $X$ is a \emph{finite} subset of $A$ if $X$ is finite and $X\subseteq A$. We will denote by $\mathrm{Sub}_{\mathrm{f}}(A)$ or, sometimes, by $\mathrm{Fin}(A)$ the set of all $S$-sorted sets $X$ in $\mathrm{Sub}(A)$ which are finite and, on some occasions, we will write $X\subseteq_{\mathrm{f}}A$ to indicate that $X\in \mathrm{Sub}_{\mathrm{f}}(A)$.
\end{definition}

\begin{remark}
For an object $A$ in the topos $\mathsf{Set}^{S}$ the following assertions are equivalent: (1) $A$ is finite, (2) $A$ is a finitary object of $\mathsf{Set}^{S}$, and (3) $A$ is a strongly finitary object of $\mathsf{Set}^{S}$ (for the notions of finitary and strongly finitary object of a category see~\cite{hs73}, Exercise 22E, on p.~155).

In $\mathsf{Set}^{S}$ there is another notion of finiteness: An $S$-sorted set $A$ is called $S$-\emph{finite} or \emph{locally finite}, abbreviated as $S\text{-}\mathrm{f}$, if and only if, for every $s\in S$, $A_{s}$ is finite.  We will denote by $\mathrm{Sub}_{S\text{-}\mathrm{f}}(A)$ the set of all $S$-sorted sets $X$ in $\mathrm{Sub}(A)$ which are $S$-finite and, on some occasions, we will write 
$X\subseteq_{S\text{-}\mathrm{f}}A$ to indicate that $X\in \mathrm{Sub}_{S\text{-}\mathrm{f}}(A)$. Although, unless $S$ is finite, this notion of finiteness is not categorial in nature, however, it plays a relevant role in the field of many-sorted algebra and in computer science. 
\end{remark}

\begin{remark}\label{RDelta}
Let $\mathrm{Fin}^{S}$ be the functor from $\mathsf{Set}^{S}$ to $\mathsf{Set}^{S}$ that sends an $S$-sorted set $A$ to the $S$-sorted set $(\mathrm{Fin}(A))_{s\in S}$, which is constantly $\mathrm{Fin}(A)$, and an $S$-sorted mapping $f\colon A\mor B$ to the $S$-sorted mapping $\mathrm{Fin}^{S}(f)$ from $(\mathrm{Fin}(A))_{s\in S}$ to $(\mathrm{Fin}(B))_{s\in S}$ that, for every $s\in S$, sends $X\in \mathrm{Fin}(A)$ to $f[X]$ in $\mathrm{Fin}(B)$, i.e., to the $S$-sorted mapping $(f[\cdot])_{s\in S}$, which is constantly $f[\cdot]$. Then the family $(\delta^{S,A})_{A\in \boldsymbol{\mathcal{U}}}$, where, for every $A\in \boldsymbol{\mathcal{U}}$,  $\delta^{S,A}$ is the $S$-sorted mapping from $A$ to $(\mathrm{Fin}(A))_{s\in S}$ that, for every $s\in S$,  sends $a\in A_{s}$ to $\delta^{S,A}_{s}(a) = \delta^{s,a}\in \mathrm{Fin}(A)$, is a natural transformation from $\mathrm{Id}^{\mathsf{Set}^{S}}$ to $\mathrm{Fin}^{S}$. This natural transformation will be used later on to define, for an $s\in S$ and an $S$-sorted set $X$, the $S$-sorted set of the variables of an $(X,s)$-term.
\end{remark}

\begin{definition}
Let $A$ be an $S$-sorted set. Then the \emph{support of}\index{many-sorted!support} $A$, denoted by $\mathrm{supp}_{S}(A)$, is the set $\{s\in S\mid A_{s}\neq \varnothing\}$.
\end{definition}

\begin{remark}
An $S$-sorted set $A$ is finite if and only if $\mathrm{supp}_{S}(A)$ is finite and, for every $s\in \mathrm{supp}_{S}(A)$, $A_{s}$ is finite.
\end{remark}

We next define the notion of equivalence relation on a many-sorted set and state the universal property of the corresponding quotient many-sorted set.

\begin{definition}\label{DEqvRel}
An $S$-\emph{sorted equivalence relation on}\index{many-sorted!equivalence relation} an $S$-sorted set $A$ is an $S$-sorted relation $\Phi$ on $A$ such that, for every $s\in S$, $\Phi_{s}$ is an equivalence relation on $A_{s}$. Sometimes, we will abbreviate the expression ``$S$-sorted equivalence relation on $A$'' to ``$S$-sorted equivalence on $A$'', or, if no confusion can arise, to ``equivalence on $A$''. We will denote by $\mathrm{Eqv}(A)$ the set of all $S$-sorted equivalences on $A$ (which is an algebraic closure system on $A\times A$), by $\mathbf{Eqv}(A)$ the algebraic lattice  $(\mathrm{Eqv}(A),\subseteq)$, by $\nabla_{A}$ the greatest element of $\mathbf{Eqv}(A)$, and by $\Delta_{A}$ the least element of $\mathbf{Eqv}(A)$. As for ordinary sets,
$\mathbf{Eqv}(A)$ is also an algebraic lattice, and we denote by
$\mathrm{Eg}_{A}$ the canonically associated algebraic closure
operator.  For $A$, we have that
$\mathrm{Eg}_{A}(\Phi)=(\mathrm{Eg}_{A_{s}}(\Phi_{s}))_{s\in S}$.

For an $S$-sorted equivalence relation $\Phi$ on $A$, the $S$-\emph{sorted quotient set of}\index{many-sorted!quotient} $A$ \emph{by} $\Phi$, denoted by  $A/\Phi$, is $(A_{s}/\Phi_{s})_{s\in S} = (\{[x]_{\Phi_{s}}\mid x\in A_{s}\})_{s\in S} (\subseteq (\mathrm{Sub}(A_{s}))_{s\in S})$, where, for every $s\in S$ and every $x\in A_{s}$, $[x]_{\Phi_{s}}$, the \emph{equivalence class of} $x$ \emph{with respect to} $\Phi_{s}$ (or, the $\Phi$-\emph{equivalence class of} $x$) is $\{y\in A_{s}\mid (x,y)\in \Phi_{s}\}$, and $\mathrm{pr}^{\Phi}\colon A\mor A/\Phi$, the \emph{canonical projection from} $A$ \emph{to} $A/\Phi$, is the $S$-sorted mapping $(\mathrm{pr}^{\Phi}_{s})_{s\in S}$, where, for every $s\in S$, $\mathrm{pr}^{\Phi}_{s}$ is the canonical projection from $A_{s}$ to $A_{s}/\Phi_{s}$ (which sends $x$ in $A_{s}$ to $\mathrm{pr}^{\Phi}_{s}(x) = [x]_{\Phi_{s}}$, the $\Phi_{s}$-equivalence class of $x$, in $A_{s}/\Phi_{s}$).  Moreover, if $\Psi$ is an $S$-sorted equivalence on an $S$-sorted set $B$ and $f$ an $S$-sorted mapping from $A$ to $B$, then the \emph{astriction of} $f$ \emph{to} $B/\Psi$, denoted by $\mathrm{ast}_{B/\Psi}(f)$, is $\mathrm{pr}^{\Psi}\circ f$, where $\mathrm{pr}^{\Psi}$ is the canonical projection of $B$ onto $B/\Psi$. Sometimes, we will abbreviate the expression ``$S$-sorted sorted quotient set of $A$ by $\Phi$'' to ``$S$-sorted quotient of $A$ by $\Phi$'', or, if no confusion can arise, to ``quotient of $A$ by $\Phi$''


\end{definition}

\begin{figure}
$$\xymatrix{
  A \ar @{{+>}}[r]^{\mathrm{pr}^{\Phi}}\ar[dr]_{f} & A/\Phi \ar[d]^{h}\\
 	                                              & B
 }
$$
\caption{The coastriction of a many-sorted mapping.}
\label{FCOAS}
\end{figure}

\begin{definition}
Let $A$ be an $S$-sorted set and $\Phi\in \mathrm{Eqv}(A)$. Then a \emph{transversal of} $A/\Phi$ \emph{in} $A$ is a subset $X$ of $A$ such that, for every $s\in S$ and every $a\in A_{s}$, $\mathrm{card}(X_{s}\cap [a]_{\Phi_{s}}) = 1$.
\end{definition}

\begin{remark}
For an $S$-sorted equivalence relation $\Phi$ on $A$, the set of all transversals of $A/\Phi$ in $A$ is isomorphic to the set of all cross-sections of $\mathrm{pr}_{\Phi}$, where an $S$-sorted mapping $f$ from $A/\Phi$ to $A$ is a cross-section of $\mathrm{pr}^{\Phi}$ if $\mathrm{pr}^{\Phi}\circ f = \mathrm{id}^{A/\Phi}$. Moreover, if $\Psi$ is another equivalence relation on $A$, $\Psi$ is a refinement of $\Phi$, i.e., $\Psi\subseteq \Phi$, and $X^{\Phi}$ is a transversal of $A/\Phi$ in $A$, then, for every $s\in S$ and every $a\in A_{s}$, there exists a unique $x\in X^{\Phi}_{s}$ such that $[a]_{\Psi_{s}}\subseteq [x]_{\Phi_{s}}$.
\end{remark}

We next define the concept of kernel of an $S$-sorted mapping.

\begin{definition}
let $f\colon A\mor B$ be an $S$-sorted mapping. Then the \emph{kernel of}\index{many-sorted!mapping!kernel} $f$, denoted by $\mathrm{Ker}(f)$, is the $S$-equivalence on $A$ defined as:
$$
\mathrm{Ker}(f)=(\mathrm{Ker}(f_{s}))_{s\in S}.
$$
\end{definition}

\begin{proposition}
Let $A$ be an $S$-sorted set, $\Phi$ an $S$-sorted equivalence on $A$, and $f\colon A\mor B$ an $S$-sorted mapping. Then the following statements are equivalent:
\begin{enumerate}
\item $\Phi\subseteq\mathrm{Ker}(f)$.
\item There exists an $S$-sorted mapping $h$ from $A/\Phi$ to $B$ such that the  diagram in Figure~\ref{FCOAS} 
commutes.
\end{enumerate}
If one of the above equivalent statements holds, then we will call $h$, which is univocally determined, the \emph{coastriction of} $f$ \emph{to} $A/\Phi$ and we denote it by $\mathrm{coast}_{A/\Phi}(f)$.
\end{proposition}

\begin{proposition}
Let $\mathsf{ClfdSet}^{S}$ be the category whose objects are the \emph{classified $S$-sorted sets}, i.e, the ordered pairs $(A,\Phi)$ where $A$ is an $S$-sorted set and $\Phi$ an $S$-sorted equivalence relation on $A$, and in which the set of morphisms from $(A,\Phi)$ to $(B,\Psi)$ is the set of all $S$-sorted mappings $f$ from $A$ to $B$ such that, for every $s\in S$ and every $(x,y)\in A^{2}_{s}$, if $(x,y)\in \Phi_{s}$, then $(f_{s}(x),f_{s}(y))\in \Psi_{s}$. Let $G$ be the functor from $\mathsf{Set}^{S}$ to $\mathsf{ClfdSet}^{S}$ whose object mapping sends $A$ to $(A,\Delta_{A})$ and whose morphism mapping sends $f\colon A\mor B$ to $f\colon (A,\Delta_{A})\mor (B,\Delta_{B})$. Then, for every classified $S$-sorted set $(A,\Phi)$, there exists a universal mapping from $(A,\Phi)$ to $G$, which is precisely the ordered pair $(A/\Phi,\mathrm{pr}^{\Phi})$ with $\mathrm{pr}^{\Phi}\colon (A,\Phi)\mor (A/\Phi,\Delta_{A/\Phi})$.
\end{proposition}

\begin{definition}
Let $f$ be a $\mathsf{ClfdSet}^{S}$-morphism from $(A,\Phi)$ to $(B,\Psi)$. Then we denote by $\mathrm{biast}_{A/\Phi,B/\Psi}(f)$, $f^{\Phi,\Psi}$, or, if no confusion can arise, $\overline{f}$ the $S$-sorted mapping $\mathrm{ast}_{B/\Psi}(\mathrm{coast}_{A/\Phi}(f))$ (which is identical to $\mathrm{coast}_{A/\Phi}(\mathrm{ast}_{B/\Psi}(f))$). We will call this $S$-sorted mapping the \emph{biastriction of} $f$ to $A/\Phi$ and $B/\Psi$. 
\end{definition}

We next state the quadrangular and triangular factorizations of an $S$-sorted mapping.

\begin{proposition}
Let $f\colon A\mor B$ be an $S$-sorted mapping. Then we have the following quadrangular and triangular factorizations of $f$:
$$\xymatrix{
A
\ar[r]^-{\mathrm{pr}^{\Ker(f)}}
\ar[d]_{f^{\mathrm{e}}} &
A/\Ker(f)
\ar[d]^{f^{\mathrm{m}}}
\ar[ld]_{f^{\mathrm{b}}}\\
\mathrm{Im}(f)
\ar[r]_{\mathrm{in}^{\mathrm{Im}(f)}} &
B
}
$$
where the involved $S$-sorted mappings are defined coordinatewise, i.e., for every $s\in S$, $(\mathrm{pr}^{\mathrm{Ker}(f)})_{s}$ is $\mathrm{pr}^{\mathrm{Ker}(f_{s})}$, the canonical projection from $A_{s}$ to $A_{s}/\Ker(f_{s})$, $f^{\mathrm{b}}_{s}$ the canonical isomorphism from $A_{s}/\mathrm{Ker}(f_{s})$ to $\mathrm{Im}(f_{s})$, $(\mathrm{in}^{\mathrm{Im}(f)})_{s}$ the
canonical embedding $\mathrm{in}^{\mathrm{Im}(f_{s})}$ of $\mathrm{Im}(f_{s})$ into $B_{s}$, $f^{\mathrm{e}}_{s}$ the corestriction of $f_{s}$ to $\mathrm{Im}(f_{s})$, and $f^{\mathrm{m}}_{s}$ the mapping which sends $[a]_{\Ker(f_{s})}$ in $A_{s}/\Ker(f_{s})$ to $f_{s}(a)$ in $B_{s}$.
\end{proposition}

%

\begin{definition}
Let $A$ be an $S$-sorted set, $X$ a subset of $A$, and $\Phi\in\mathrm{Eqv}(A)$. Then the $\Phi$-\emph{saturation of} $X$ (or, the \emph{saturation of} $X$ \emph{with respect to} $\Phi$), denoted by $[X]^{\Phi}$, is the $S$-sorted set defined, for every $s\in S$, as follows:
$$
\textstyle[X]^{\Phi}_{s} = \{a\in A_{s}\mid
X_{s}\cap [a]_{\Phi_{s}}\neq\varnothing\}= \bigcup_{x\in X_{s}}[x]_{\Phi_{s}} = [X_{s}]^{\Phi_{s}}.
$$
Let $X$ be a subset of $A$ and $\Phi\in\mathrm{Eqv}(A)$. Then we will say that $X$ is $\Phi$-\emph{saturated} if and only if $X = [X]^{\Phi}$. We will denote by $\Phi\text{-}\mathrm{Sat}(A)$ the subset of $\mathrm{Sub}(A)$ defined as $\Phi\text{-}\mathrm{Sat}(A) = \{X\in \mathrm{Sub}(A)\mid X = [X]^{\Phi}\}$.
\end{definition}

\begin{remark}
Let $A$ be an $S$-sorted set and $\Phi\in\mathrm{Eqv}(A)$. Then, for a subset $X$ of
$A$, we have that the $\Phi$-saturation of $X$ is $(\mathrm{pr}^{\Phi})^{-1}[\mathrm{pr}^{\Phi}[X]]$. Therefore, $X$ is $\Phi$-saturated if and only if $X \supseteq [X]^{\Phi}$. Besides, $X$ is $\Phi$-saturated if and only if there exists a $\mathcal{Y}\subseteq A/\Phi$ such that $X = (\mathrm{pr}^{\Phi})^{-1}[\mathcal{Y}]$.
\end{remark}

\begin{proposition}\label{PSatIncl}
Let $A$ be an $S$-sorted set and $\Phi$, $\Psi\in\mathrm{Eqv}(A)$. Then
$$
\Phi\subseteq \Psi\, \text{if and only if }\, \forall X\subseteq A\;([[X]^{\Psi}]^{\Phi}=[X]^{\Psi}).
$$
Moreover, for a sort $s\in S$, we have that
$$
\Phi_{s}\subseteq \Psi_{s}\, \text{if and only if }\, \forall X\subseteq A_{s}\;([[X]^{\Psi_{s}}]^{\Phi_{s}}=[X]^{\Psi_{s}}).
$$
\end{proposition}

\begin{proof}
We restrict ourselves to proving the first assertion. Let us assume that $\Phi\subseteq \Psi$ and let $X$ be a subset of $A$. In order to prove that $[[X]^{\Psi}]^{\Phi}=[X]^{\Psi}$ it suffices to verify that $[[X]^{\Psi}]^{\Phi}\subseteq [X]^{\Psi}$. Let $s$ be an element of $S$. Then, by definition, $a\in [[X]^{\Psi}]^{\Phi}_{s}$ if and only if there exists some $b\in [X]^{\Psi}_{s}$ such that $a\in [b]_{\Phi_{s}}$. Since $\Phi\subseteq \Psi$, we have that $a\in [b]_{\Psi_{s}}$, therefore $a\in [X]^{\Psi}_{s}$.

To prove the converse, let us assume that $\Phi\not\subseteq \Psi$. Then there exists some sort $s\in S$ and elements $a$, $b$ in $A_{s}$ such that $(a,b)\in\Phi_{s}$ and $(a,b)\not\in \Psi_{s}$. Hence $b$ does not belong to $[\delta^{s,[a]_{\Psi_{s}}}]^{\Psi}_s$, whereas it does belong to $[[\delta^{s,[a]_{\Psi_{s}}}]^{\Psi}]^{\Phi}_{s}$. It follows that $[\delta^{s,[a]_{\Psi_{s}}}]^{\Psi}\neq [[\delta^{s,[a]_{\Psi_{s}}}]^{\Psi}]^{\Phi}$.
\end{proof}

\begin{corollary}\label{CSatIncl}
Let $A$ be an $S$-sorted set and $\Phi$, $\Psi\in\mathrm{Eqv}(A)$. If $\Phi\subseteq \Psi$, then  $\Psi\text{-}\mathrm{Sat}(A)\subseteq\Phi\text{-}\mathrm{Sat}(A)$. Moreover, for $s\in S$ and $L\subseteq A_{s}$, if $\Phi_{s}\subseteq \Psi_{s}$ and $L = [L]^{\Psi_{s}}$, then  $L = [L]^{\Phi_{s}}$.
\end{corollary}

\begin{remark}
If, for an $S$-sorted set $A$, we denote by $(\cdot)\text{-}\mathrm{Sat}(A)$ the mapping from $\mathrm{Eqv}(A)$ to $\mathrm{Sub}(\mathrm{Sub}(A))$ which sends $\Phi$ to $\Phi\text{-}\mathrm{Sat}(A)$, then the above corollary means that $(\cdot)\text{-}\mathrm{Sat}(A)$ is an antitone ($\equiv$ order-reversing) mapping from the ordered set $(\mathrm{Eqv}(A),\subseteq)$ to the ordered set $(\mathrm{Sub}(\mathrm{Sub}(A)),\subseteq)$.
\end{remark}

\begin{proposition}\label{PSatNabla}
Let $A$ be an $S$-sorted set and $X\subseteq A$. Then $X\in \nabla_{A}\text{-}\mathrm{Sat}(A)$ if and only if, for every $s\in S$, if $s\in \mathrm{supp}_{S}(X)$, then $X_{s} = A_{s}$.
\end{proposition}

\begin{proof}
Let us suppose that there exists a $t\in S$ such that $X_{t}\neq \varnothing$ and $X_{t}\neq A_{t}$. Then, since $[X]_{t}^{\nabla_{A}} = \bigcup_{x\in X_{t}}[x]_{\nabla_{A_{t}}}$ and $X_{t}\neq \varnothing$, we have that, for some $y\in X_{t}$, $[y]_{\nabla_{A_{t}}} = A_{t}$. But $X_{t}\subset A_{t}$. Hence $[X]_{t}^{\nabla_{A}}\neq X_{t}$. Therefore $X\not\in \nabla_{A}\text{-}\mathrm{Sat}(A)$.

The converse implication is straightforward.
\end{proof}

\begin{remark}
Let $A$ be an $S$-sorted set. Then, from the above proposition, it follows that $\varnothing^{S}$, $A\in \nabla_{A}\text{-}\mathrm{Sat}(A)$. Moreover, for every subset $T$ of $S$, we have that $\bigcup_{t\in T}\delta^{t,A_{t}}\in \nabla_{A}\text{-}\mathrm{Sat}(A)$.
\end{remark}

\begin{proposition}
Let $A$ be an $S$-sorted set, $X\subseteq A$, and $\Phi$, $\Psi\in\mathrm{Eqv}(A)$. Then
$[X]^{\Phi\cap\Psi}\subseteq[X]^{\Phi}\cap[X]^{\Psi}$.
\end{proposition}

\begin{proof}
Let $s$ be a sort in $S$ and $b\in[X]^{\Phi\cap \Psi}_{s}$. Then, by definition, there exists an $a\in X_{s}$ such that $(a,b)\in (\Phi\cap\Psi)_{s} = \Phi_{s}\cap\Psi_{s}$. Hence, $(a,b)\in \Phi_{s}$ and $(a,b)\in \Psi_{s}$. Therefore $b\in [X]^{\Phi}_{s}$ and $b\in[X]^{\Psi}_{s}$. Consequently, $b\in ([X]^{\Phi}\cap [X]^{\Psi})_{s}$. Thus $[X]^{\Phi\cap\Psi}\subseteq[X]^{\Phi}\cap[X]^{\Psi}$.
\end{proof}

\begin{corollary}
Let $A$ be an $S$-sorted set and $\Phi$, $\Psi\in\mathrm{Eqv}(A)$. Then we have that  $\Phi\text{-}\mathrm{Sat}(A)\cap\Psi\text{-}\mathrm{Sat}(A)\subseteq(\Phi\cap\Psi)\text{-}\mathrm{Sat}(A)$.
\end{corollary}

We next state that the set $\Phi\text{-}\mathrm{Sat}(A)$ is the set of all fixed points of a suitable operator on $A$, i.e., of an endomapping of $\mathrm{Sub}(A)$.

\begin{proposition}\label{PSatOp}
Let $A$ be an $S$-sorted set and $\Phi\in\mathrm{Eqv}(A)$. Then the mapping $[\cdot]^{\Phi}$ from $\mathrm{Sub}(A)$ to $\mathrm{Sub}(A)$ that sends $X$ in $\mathrm{Sub}(A)$ to $[\cdot]^{\Phi}(X) = [X]^{\Phi}$ in $\mathrm{Sub}(A)$
is a completely additive closure operator on $A$. Moreover, for every nonempty set $I$ in $\boldsymbol{\mathcal{U}}$ and every $I$-indexed family $(X^{i})_{i\in I}$ in $\mathrm{Sub}(A)$, $[\bigcap_{i\in I}X^{i}]^{\Phi} \subseteq \bigcap_{i\in I}[X^{i}]^{\Phi}$ (and, obviously, $[A]^{\Phi} = A$), and, for every $X\subseteq A$, if $X = [X]^{\Phi}$, then  $\complement_{A}X = [\complement_{A}X]^{\Phi}$. Besides, $[\cdot]^{\Phi}$ is uniform, i.e., is such that, for every $X$, $Y\subseteq A$, if  $\mathrm{supp}_{S}(X) = \mathrm{supp}_{S}(Y)$, then $\mathrm{supp}_{S}([X]^{\Phi}) = \mathrm{supp}_{S}([Y]^{\Phi})$---hence, in particular, $[\cdot]^{\Phi}$ is a uniform algebraic closure operator on $A$. And $\Phi\text{-}\mathrm{Sat}(A) = \mathrm{Fix}([\cdot]^{\Phi})$, where $\mathrm{Fix}([\cdot]^{\Phi})$ is the set of all fixed point of the operator $[\cdot]^{\Phi}$.
\end{proposition}


\begin{proposition}\label{PSatCABA} Let $A$ be an $S$-sorted set and $\Phi\in\mathrm{Eqv}(A)$. Then the ordered pair
$\Phi\text{-}\mathbf{Sat}(A) = (\Phi\text{-}\mathrm{Sat}(A),\subseteq)$ is a complete atomic Boolean algebra.
\end{proposition}

\begin{proof}
The proof is straightforward and we leave it to the reader. We only point out that the atoms of $\Phi\text{-}\mathbf{Sat}(A)$ are precisely the deltas of Kronecker $\delta^{t, [x]_{\Phi_t}}$, for some $t\in S$ and some $x\in A_{t}$, and that, obviously, every $\Phi$-saturated subset $X$ of $A$ is the join ($\equiv$ union) of all atoms smaller than $X$.
\end{proof}

\begin{proposition}\label{POptLift}
Let $A$ be an $S$-sorted set, $(B,\Psi)$ a classified $S$-sorted set and $f\colon A\mor B$ an $S$-sorted mapping; situation that we will indicate by:
$$
  f\colon A\mor (B,S).
$$
Then there exists an \emph{optimal lift of} $\Psi$ \emph{through} $f$, i.e., there exists an $S$-sorted equivalence on $A$, denoted by $\mathrm{L}^{f}(\Psi)$, the \emph{optimal lift of} $\Psi$ \emph{through} $f$, such that 
\begin{enumerate}
\item $f\colon (A,\mathrm{L}^{f}(\Psi))\mor (B,\Psi)$ is a morphism and,
\item for every  classified $S$-sorted set $(D,\Upsilon)$ and every $S$-sorted mapping $g\colon D\mor A$, 
if $f\circ g\colon (D,\Upsilon)\mor (B,\Psi)$ is a morphism, then
$g\colon (D,\Upsilon)\mor(A,\mathrm{L}^{f}(\Psi))$ is a morphism.
\end{enumerate}
Moreover, we have that:
\begin{enumerate}
\item For every $S$-sorted equivalence $\Phi\in \mathrm{Eqv}(A)$:
$$
  \mathrm{L}^{\mathrm{id}^{A}}(\Phi) = \Phi.
$$
\item If $f\colon A\mor B$, $g\colon B\mor C$ are $S$-sorted mappings and $\Omega\in
\mathrm{Eqv}(C)$, then:
$$
  \mathrm{L}^{g\circ f}(\Omega) = \mathrm{L}^{f}(\mathrm{L}^{g}(\Omega)).
$$
\end{enumerate}
\end{proposition}


\begin{remark}\label{ROptLiftDiag}
Let $A$ be an $S$-sorted set, $(B,\Psi)$ a classified $S$-sorted set and $f\colon A\mor B$ an $S$-sorted mapping. Then there exists a unique $S$-sorted mapping $f^{\mathrm{L}^{f}(\Psi),\Psi}$ from $A/\mathrm{L}^{f}(\Psi)$ to $B/\Psi$ such that $\mathrm{pr}^{\Psi}\circ f = f^{\mathrm{L}^{f}(\Psi),\Psi}\circ \mathrm{pr}^{\mathrm{L}^{f}(\Psi)}$, i.e., the following diagram commutes
$$\xymatrix{
A
\ar[r]^-{f}
\ar[d]_{\mathrm{pr}^{\mathrm{L}^{f}(\Psi)}} &
B
\ar[d]^{\mathrm{pr}^{\Psi}}
\\
A/\mathrm{L}^{f}(\Psi)
\ar[r]_{f^{\mathrm{L}^{f}(\Psi),\Psi}} &
B/\Psi
}
$$
\end{remark}

\begin{proposition}\label{PClfd}
The functor from $\mathsf{ClfdSet}$ to $\mathsf{ClfdSet}^{S}$, whose object mapping assigns to a classified set $(X,\Theta)$ the classified $S$-sorted set $((X)_{s\in S},(\Theta)_{s\in S})$ and whose morphism mapping sends a morphism $f$ from $(X,\Theta)$ to $(Y,\Omega)$ to the morphism $(f)_{s\in S}$ from $((X)_{s\in S},(\Theta)_{s\in S})$ to $((Y)_{s\in S},(\Omega)_{s\in S})$, has a left adjoint.
\end{proposition}

\begin{proof}
Let $(A,\Phi)$ be a classified $S$-sorted set. Then $(\coprod A, \Phi^{\coprod A})$, where $\Phi^{\coprod A}$ is the equivalence relation on $\coprod A$ defined as
$$
\textstyle
\Phi^{\coprod A} = \{((x,s),(y,t))\in ( \coprod A)^{2}\mid s = t \And (x,y)\in \Phi_ {s}\},
$$
together with the canonical embedding of $(A,\Phi)$ into $((\coprod A)_{s\in S},(\Phi^{\coprod A})_{s\in S})$ is a universal morphism from $(A,\Phi)$ to the  functor from $\mathsf{ClfdSet}$ to $\mathsf{ClfdSet}^{S}$. The details are left to the reader.
\end{proof}

\section{Free monoids}

We next define the concept of free monoid on a set and several notions associated with it that will be used afterwards to construct the free algebra on an $S$-sorted set and to define various substitution operators.

\begin{definition}
Let $A$ be a set. The \emph{free monoid on} $A$, denoted by $\mathbf{A}^{\star}$, is $(A^{\star},\curlywedge,\lambda)$, where $A^{\star}$, the set of all \emph{words on} $A$, is $\bigcup_{n\in\mathbb{N}}\mathrm{Hom}(n,A)$, the set of all mappings $\mathbf{a}\colon n\mor A$ from some $n\in \mathbb{N}$ to $A$. A word $\mathbf{a}\in A^{\star}$ is usually denoted as a sequence $(a_{i})_{i\in\bb{\mathbf{a}}}$, where, for $i\in\bb{\mathbf{a}}$, $a_{i}$ is the letter in $A$ satisfying $\mathbf{a}(i)=a_{i}$. Furthermore, $\curlywedge$, the \emph{concatenation} of words on $A$, is the binary operation on $A^{\star}$ which sends a pair of words $(\mathbf{a},\mathbf{b})$ on $A$ to the mapping $\mathbf{a}\curlywedge \mathbf{b}$ from $\bb{\mathbf{a}}+\bb{\mathbf{b}}$ to $A$, where $\bb{\mathbf{a}}$ and $\bb{\mathbf{b}}$ are the lengths ($\equiv$ domains) of the mappings $\mathbf{a}$ and $\mathbf{b}$, respectively, defined as follows:
$$
\mathbf{a}\bconcat \mathbf{b}
\nfunction
{\bb{\mathbf{a}}+\bb{\mathbf{b}}}{A}
{i}{
\begin{cases}
  a_{i}, & \text{if $0\leq i < \bb{\mathbf{a}}$;}\\
  b_{i-\bb{\mathbf{a}}}, & \text{if $\bb{\mathbf{a}}\leq i < \bb{\mathbf{a}}+\bb{\mathbf{b}}$,}
\end{cases}
   }
$$
and $\lambda$, the \emph{empty word on} $A$, is the unique mapping from $\varnothing$ to $A$. We will denote by $\eta^{A}$ the mapping from $A$ to $A^{\star}$ that sends $a\in A$ to $(a)\in A^{\star}$, i.e., to the mapping $(a)\colon 1\mor A$ that sends $0$ to $a$. The ordered pair $(\mathbf{A}^{\star},\eta^{A})$ is a universal morphism from $A$ to the forgetful functor from the category $\mathsf{Mon}$, of monoids, to $\mathsf{Set}$.
\end{definition}

\begin{remark}
For a word $\mathbf{a}\in A^{\star}$, $\bb{\mathbf{a}}$, the length of $\mathbf{a}$, is the value at $\mathbf{a}$ of the unique homomorphism $\bb{\,\cdot\,}$ from $\mathbf{A}^{\star}$ to $(\mathbb{N},+,0)$, the additive monoid of the natural numbers, such that $\bb{\,\cdot\,}\circ \eta^{A} = \kappa_{1}$, where $\kappa_{1}$ is the mapping from $A$ to $\mathbb{N}$ constantly $1$. Note that, for every $n\in \mathbb{N}$, $\bb{\,\cdot\,}$ sends $\mathbf{a}\in \mathrm{Hom}(n,A)$ to $n$. Thus, for the family of mappings $(\kappa_{n})_{n\in \mathbb{N}}$, where, for every $n\in \mathbb{N}$, $\kappa_{n}$ is the mapping from $\mathrm{Hom}(n,A)$ to $\mathbb{N}$ constantly $n$, and by applying the universal property of the coproduct, we have $\bb{\,\cdot\,} = [\kappa_{n}]_{n\in \mathbb{N}}$.
\end{remark}

%

\begin{definition}\label{DSubw}
Let $\mathbf{a}$ and $\mathbf{a}'$ be words in $A^{\star}$ and $k$, $l\in\bb{\mathbf{a}}$ such that $k\leq l$. We will say that $\mathbf{a}'$ is a \emph{subword} of $\mathbf{a}$ \emph{beginning at position} $k$ and \emph{ending at position} $l$ if there are words $\mathbf{b}$ and $\mathbf{c}$ such that $\mathbf{a} = \mathbf{b}\bconcat \mathbf{a}'\bconcat \mathbf{c}$, $\bb{\mathbf{b}} = k$, and $\bb{\mathbf{c}} = (\bb{\mathbf{a}}-1)-l$. For $\mathbf{a}\in A^{\star}$ and $k$, $l\in\bb{\mathbf{a}}$ such that $k\leq l$, we will denote by $\mathbf{a}^{k,l}$ the word $\mathbf{a}'$  beginning at position $k$ and ending at position $l$. In particular, we will denote by $\mathbf{a}^{0,l}$ the word $\mathbf{a}'$  beginning at position $0$ and ending at position $l$, which is an \emph{initial segment} of $\mathbf{a}$, and by $\mathbf{a}^{k,\bb{\mathbf{a}}-1}$ the word $\mathbf{a}'$  beginning at position $k$ and ending at position $\bb{\mathbf{a}}-1$, which is a \emph{final segment} of $\mathbf{a}$. We will call the word $\mathbf{a}^{0,\bb{\mathbf{a}}-2}$ the \emph{maximal prefix} of $\mathbf{a}$. Moreover, we will call $\mathbf{a}^{0,0}$, which is, essentially, the first letter of $\mathbf{a}$, the \emph{head} of $\mathbf{a}$ and $\mathbf{a}^{\bb{\mathbf{a}}-1,\bb{\mathbf{a}}-1}$, which is, essentially, the last letter of $\mathbf{a}$, the \emph{tail} of $\mathbf{a}$. A word $\mathbf{a}$ may have several subwords equal to $\mathbf{a}'$. In that case, the equation $\mathbf{a} = \mathbf{b}\bconcat \mathbf{a}'\bconcat \mathbf{c}$ has several solutions $(\mathbf{b},\mathbf{c})$.
If the pairs $(\mathbf{b}_{i},\mathbf{c}_{i})$ $(i\in n)$ are all solutions of $\mathbf{a} = \mathbf{b}\bconcat \mathbf{a}'\bconcat \mathbf{c}$ and if $\bb{\mathbf{b}_{0}}< \bb{\mathbf{b}_{1}}<\cdots<\bb{\mathbf{b}_{n-1}}$, then $(\mathbf{b}_{i},\mathbf{c}_{i})$ determines the $i$-th occurrence of $\mathbf{a}'$ in $\mathbf{a}$. The solution $(\mathbf{b},\mathbf{c})$ in which either $\mathbf{b}$ or $\mathbf{c}$ is $\lambda$ is not excluded.
Let $\mathbf{a}$ be a word in $A^{\star}$ and $a\in A$. We will say that \emph{$a$ occurs in} $\mathbf{a}$ if there are words $\mathbf{b}$, $\mathbf{c}$ in $A^{\star}$ such that $\mathbf{a} = \mathbf{b}\bconcat(a)\bconcat \mathbf{c}$. Note that $a$ occurs in $\mathbf{a}$ if and only if there exists an $i\in \bb{\mathbf{a}}$ such that $\mathbf{a}(i) = a$. We will denote by $\bb{\mathbf{a}}_{a}$ the natural number $\mathrm{card}(\{i\in \bb{\mathbf{a}}\mid \mathbf{a}(i) = a\}) = \mathrm{card}(\mathbf{a}^{-1}[\{a\}])$, i.e., the number of occurrences of $a$ in $\mathbf{a}$. Moreover, we let $(i_{j})_{j\in\bb{\mathbf{a}}_{a}}$ stand for the  enumeration in ascending order of the occurrences of $a$ in $\mathbf{a}$. Thus $(i_{j})_{j\in\bb{\mathbf{a}}_{a}}$ is the order embedding of $(\bb{\mathbf{a}}_{a},<)$ into $(\bb{\mathbf{a}},<)$ defined  recursively as follows:
$$i_{0} = \mathrm{min}\{i\in \bb{\mathbf{a}}\mid \mathbf{a}(i) = a\},$$
and, for $j\in [1,\bb{\mathbf{a}}_{a}-1]$, 
$$
i_{j} = \mathrm{min}\{i\in \bb{\mathbf{a}}-\{i_{0},\ldots,i_{j-1}\}\mid \mathbf{a}(i) = a\}.
$$
If, for every $j\in \bb{\mathbf{a}}_{a}$,  the pairs $(\mathbf{b}_{i_{j}},\mathbf{c}_{i_{j}})$  are all solutions of $\mathbf{a} = \mathbf{b}_{i_{j}}\bconcat (a)\bconcat \mathbf{c}_{i_{j}}$ and if $\bb{\mathbf{b}_{i_{0}}}< \bb{\mathbf{b}_{i_{1}}}<\cdots<\bb{\mathbf{b}_{i_{\bb{\mathbf{a}}_{a}-1}}}$, then $(\mathbf{b}_{i_{j}},\mathbf{c}_{i_{j}})$ determines the $j$-th occurrence of $(a)$ in $\mathbf{a}$ and we will say that \emph{$a$ occurs at the $i_{j}$-th place of $\mathbf{a}$}.
\end{definition}


\begin{remark}
Let $a$ be an element of $A$. Then, for $\mathbf{a}\in A^{\star}$, $\bb{\mathbf{a}}_{a}$, the number of occurrences of $a$ in $\mathbf{a}$, is the value at $\mathbf{a}$ of the unique homomorphism $\bb{\,\cdot\,}_{a}$ from $\mathbf{A}^{\star}$ to $(\mathbb{N},+,0)$ such that $\bb{\,\cdot\,}_{a}\circ\eta^{A} = \delta_{a}$, where $\delta_{a}$ is the mapping from $A$ to $\mathbb{N}$ that sends $a\in A$ to $1\in\mathbb{N}$ and $b\in A-\{a\}$ to $0\in \mathbb{N}$. Therefore, since, for every $\mathbf{a}\in A^{\star}$, the $A$-indexed family $(\bb{\mathbf{a}}_{a})_{a\in A}$ in $\mathbb{N}$ is such that $\bb{\mathbf{a}}_{a} = 0$ for all but a finite number of elements $a$ in $A$, i.e., is such that $\mathrm{card}(\{a\in A\mid \bb{\mathbf{a}}_{a}\neq 0\})<\aleph_{0}$, we have that $\sum_{a\in A}\bb{\mathbf{a}}_{a}\in \mathbb{N}$ and, obviously, for every $\mathbf{a}\in A^{\star}$, $\bb{\mathbf{a}} = \sum_{a\in A}\bb{\mathbf{a}}_{a}$, i.e., $\bb{\,\cdot\,} = \sum_{a\in A}\bb{\,\cdot\,}_{a}$.
\end{remark}

\section{Artinian ordered sets}

In this subsection we recall the notion of Artinian ordered set.

\begin{definition}\label{DPosArt}
A \emph{preordered} set is an ordered pair $\mathbf{A}=(A,\leq)$ in which $A$ is a set and $\leq$ a reflexive and transitive binary relation on $A$. An \emph{ordered} set is an ordered pair $\mathbf{A}=(A,\leq)$ in which $A$ is a set and $\leq$ a reflexive, antisymmetric and transitive binary relation on $A$. 
A \emph{strictly ordered} set is an ordered pair $\mathbf{A}=(A,<)$ in which $A$ is a set and $<$ an irreflexive and transitive binary relation on $A$ (recall that, for every set $A$, there exists a  biunivocal correspondence between the set of the orderings on $A$ and the set of the strict orderings on $A$). We will say that an ordered set $\mathbf{A}=(A,\leq)$ is \emph{Artinian} if in $\mathbf{A}$ there is not any strictly decreasing $\omega_{0}$-chain or, what is equivalent, if every nonempty subset of $A$ has a minimal element. Finally, a subset $X$ of $A$ is a \emph{lower subset} of $\mathbf{A}$ if it is downward closed, i.e., if, for every $x\in X$ and every $a\in A$, if $a\leq x$, then $a\in X$.
\end{definition}

\begin{remark}
An ordered set $\mathbf{A}$ is Artinian if and only if, for every $\omega_{0}$-indexed family $(a_{n})_{n\in\omega_{0}}$ in $A$, if, for every $n\in\omega_{0}$, $a_{n+1}\leq a_{n}$, then there exists an $m\in\omega_{0}$ such that, for every $p\in \omega_{0}$, $a_{m} = a_{m+p}$. In other words, an ordered set $\mathbf{A}$ is Artinian if and only if every decreasing $\omega_{0}$-chain in $\mathbf{A}$ becomes stationary. Moreover, an ordered set $\mathbf{A}$ is Artinian if and only if, for every $X\subseteq A$, if $\mathrm{Min}(\mathbf{A})$, the set of all minimals of $\mathbf{A}$, is included in $X$ and, for every $a\in A$, if from $\{b\in A\mid b<a\}\subseteq X$ it follows that $a\in X$, then $X=A$. This last condition enables us to carry out not only proofs by Artinian induction but also definitions by Artinian recursion  (see~\cite{agk67} for details). The latter will be fundamental to this work, specifically when we define the $S$-sorted Curry-Howard mappings.
\end{remark}

\begin{remark}
Let us recall that in Set Theory (see, e.g.,~\cite{end77}) a binary relation $\Phi$ on a set $A$ is \emph{well-founded} if every nonempty subset $X$ of $A$ has a $\Phi$-minimal element, i.e., there exists an $m\in X$ such that, for every $x\in X$, $(x,m)\nin \Phi$. If the binary relation $\Phi$ on $A$ is well-founded, then it is irreflexive, i.e., for every $a\in A$, $(a,a)\nin \Phi$ (otherwise, $\Phi$ would not be well-founded). On the other hand, the binary relation $\Phi$ on $A$ is well-founded if and only if there is no $\omega_{0}$-family $(a_{n})_{n\in \mathbb{N}}$ in $A$ such that, for every $n\in \mathbb{N}$, $(a_{n+1},a_{n})\in \Phi$. Finally, if the binary relation $\Phi$ on $A$ is well-founded, then $\Phi^{+}$, its transitive closure---which is $\bigcup_{n\in \mathbb{N}-1}\Phi^{n}$---, is also well-founded and, consequently, $\Phi^{+}$ is irreflexive and transitive on $A$ (i.e., $\Phi^{+}$ is an strict order on $A$) and $(A,\Phi)$ is an Artinian ordered set.
\end{remark}




\section{Many-sorted algebras}

Our next aim is to provide those notions from the field of many-sorted universal algebra that will be used afterwards.
We will specially focus on the constructive description of the subalgebra generating many-sorted operator, the congruence generating many-sorted operator, the many-sorted translations, the (absolutely) free many-sorted algebra on a many-sorted set, and the substitution operators associated to a free many-sorted algebra.

\begin{convention} In what follows, for a set of sorts $S$, an arbitrary word on $S^{\star}$ will be denoted by $\mathbf{s}$, i.e., a lower case bold type $s$. The letter $s$ will be used to represent an arbitrary letter in $S$.
\end{convention}

\begin{definition}\label{DSig}
An $S$-\emph{sorted signature} is a mapping $\Sigma$ from $S^{\star}\times S$ to $\boldsymbol{\mathcal{U}}$ which sends a pair $(\mathbf{s},s)\in S^{\star}\times S$ to the set $\Sigma_{\mathbf{s},s}$ of the \emph{formal operations} of \emph{arity} $\mathbf{s}$, \emph{sort} (or \emph{coarity}) $s$, and \emph{rank} (or \emph{biarity}) $(\mathbf{s},s)$.

\end{definition}



\begin{assumption}
From now on $\Sigma$ stands for an $S$-sorted signature, fixed once and for all.
\end{assumption}

%

We shall now give precise definitions of the concepts of many-sorted algebra and homomorphism between many-sorted algebras.

\begin{definition}\label{DAlg}
The $S^{\star}\times S$-sorted set of the \emph{finitary operations on} an $S$-sorted set $A$ is $(\mathrm{Hom}(A_{\mathbf{s}},A_{s}))_{(\mathbf{s},s)\in S^{\star}\times S}$, where, for every $\mathbf{s}\in S^{\star}$, $A_{\mathbf{s}} = \prod_{j\in \bb{\mathbf{s}}}A_{s_{j}}$, with $\bb{\mathbf{s}}$ denoting the length of the word $\mathbf{s}$ (if $\mathbf{s} = \lambda$, then $A_{\lambda}$ is a final set). Sometimes we let $\mathrm{O}_{\mathrm{H}_{S}}(A)$ stand for 
$(\mathrm{Hom}(A_{\mathbf{s}},A_{s}))_{(\mathbf{s},s)\in S^{\star}\times S}$. A \emph{structure of} $\Sigma$-\emph{algebra on} an $S$-sorted set  $A$ is a family $(F_{\mathbf{s},s})_{(\mathbf{s},s)\in S^{\star}\times S}$, denoted by $F$, where, for $(\mathbf{s},s)\in S^{\star}\times S$, $F_{\mathbf{s},s}$ is a mapping from $\Sigma_{\mathbf{s},s}$ to $\mathrm{Hom}(A_{\mathbf{s}},A_{s})$ (if $(\mathbf{s},s) = (\lambda,s)$ and $\sigma\in \Sigma_{\lambda,s}$, then $F_{\lambda,s}(\sigma)$ picks out an element of $A_{s}$).  A \emph{many-sorted} $\Sigma$-\emph{algebra} (or, to abbreviate, $\Sigma$-\emph{algebra}) is a pair $(A,F)$, denoted by $\mathbf{A}$, where $A$ is an $S$-sorted set and $F$ a structure of $\Sigma$-algebra on $A$. For a pair $(\mathbf{s},s)\in S^{\star}\times S$ and a formal operation $\sigma\in \Sigma_{\mathbf{s},s}$, in order to simplify the notation, the operation $F_{\mathbf{s},s}(\sigma)$ from $A_{\mathbf{s}}$ to $A_{s}$ will be written as $F_{\sigma}$. In some cases, to avoid mistakes, we will denote by $F^{\mathbf{A}}$ the structure of $\Sigma$-algebra on $A$, and, for $(\mathbf{s},s)\in S^{\star}\times S$ and $\sigma\in \Sigma_{\mathbf{s},s}$, by $F^{\mathbf{A}}_{\sigma}$, or simply by $\sigma^{\mathbf{A}}$, the corresponding operation. Moreover, for $s\in S$ and $\sigma\in\Sigma_{\lambda,s}$, we will, usually, denote by $\sigma^{\mathbf{A}}$ the value of the mapping $F^{\mathbf{A}}_{\sigma}\colon A_{\lambda}\mor A_{s}$ at the unique element in $A_{\lambda}$.

A $\Sigma$-\emph{homomorphism} (or, to abbreviate, \emph{homomorphism}) from $\mathbf{A}$ to $\mathbf{B}$, where $\mathbf{B} = (B,G)$, is a triple $(\mathbf{A},f,\mathbf{B})$, denoted by $f\colon \mathbf{A}\mor \mathbf{B}$, where $f$ is an $S$-sorted mapping from $A$ to $B$ such that, for every $(\mathbf{s},s)\in S^{\star}\times S$, every  $\sigma\in \Sigma_{\mathbf{s},s}$, and every $(a_{j})_{j\in \bb{\mathbf{s}}}\in A_{\mathbf{s}}$, we have 
$
f_{s}(\sigma^{\mathbf{A}}((a_{j})_{j\in \bb{\mathbf{s}}})) = \sigma^{\mathbf{B}}(f_{\mathbf{s}}((a_{j})_{j\in \bb{\mathbf{s}}})),
$
where $f_{\mathbf{s}}$ is the mapping $\prod_{j\in \bb{\mathbf{s}}}f_{s_{j}}$ from $A_{\mathbf{s}}$ to $B_{\mathbf{s}}$ that sends $(a_{j})_{j\in \bb{\mathbf{s}}}$ in $A_{\mathbf{s}}$ to $(f_{s_{j}}(a_{j}))_{j\in \bb{\mathbf{s}}}$ in $B_{\mathbf{s}}$. For a $\Sigma$-algebra $\mathbf{A}$ the triple $(\mathbf{A},\mathrm{id}^{A},\mathbf{A})$, where $\mathrm{id}^{A}$ is the identity at $A$, which is a homomorphism from $\mathbf{A}$ to $\mathbf{A}$, will be denoted by $\mathrm{id}^{\mathbf{A}}$ and will be called the \emph{identity} homomorphism at $\mathbf{A}$. We will denote by $\mathsf{Alg}(\Sigma)$ the category of $\Sigma$-algebras and homomorphisms and by $\mathrm{Alg}(\Sigma)$ the set of objects  of $\mathsf{Alg}(\Sigma)$.

We will denote by $\mathbf{1}^{S}$ or, to abbreviate, by $\mathbf{1}$, the (standard) final $\Sigma$-algebra.
\end{definition}

We add some further results concerning direct products and homomorphisms, which will be of great importance later on.

\begin{definition}\label{DAlgProd}
Let $(\mathbf{A}^{i})_{i\in I}$ be a family of $\Sigma$-algebras, where we agree that, for every $i\in I$, $\mathbf{A}^{i}=(A^{i},F^{i})$. Then
\begin{enumerate}
\item The \emph{product}
      of $(\mathbf{A}^{i})_{i\in I}$, $\prod_{i\in I}\mathbf{A}^{i}$, is
      the  $\Sigma$-algebra which has as $S$-sorted underlying set $\prod_{i\in I}A^{i}$, 
      and where, for every $(\mathbf{s},s)\in S^{\star}\times S$ and every $\sigma\in \Sigma_{\mathbf{s},s}$,
      the structural operation $\sigma^{\prod_{i\in I}\mathbf{A}^{i}}$ is defined as follows:
      $$
      \sigma^{\prod_{i\in I}\mathbf{A}^{i}}
      \nfunction
      { (\prod_{i\in I}A^{i})_{\mathbf{s}}}
      {\prod_{i\in I}A^{i}_{s}}
      {(a_{j})_{j\in \bb{\mathbf{s}}}}
      {(\sigma^{\mathbf{A}_{i}}((a_{j}(i))_{j\in \bb{\mathbf{s}}})_{i\in I}}.
      $$
\item For every $i\in I$, the $i$-th \emph{canonical projection}
      is the homomorphism from $\prod_{i\in I}\mathbf{A}^{i}$ to $\mathbf{A}^{i}$
      determined by the $S$-sorted mapping $\pr^{i}$ which, for every $s\in S$, is
      defined as follows
      $$
      \pr^{i}_{s}
      \nfunction
      {\prod_{i\in I}A^{i}_{s}}
      {A^{i}_{s}}
      {(a_{i})_{i\in I}}
      {a_{i}}
      $$
\end{enumerate}
\end{definition}

\begin{proposition}\label{PAlgProd}
Let $(\mathbf{A}^{i})_{i\in I}$ be a family of $\Sigma$-algebras.
Then the pair $(\prod_{i\in I}\mathbf{A}^{i},(\pr^{i})_{i\in I})$
is a product in $\mathsf{Alg}(\Sigma)$.
\end{proposition}

We next introduce the support of a many-sorted algebra and the definition of finiteness of a many-sorted algebra.

\begin{definition}
Let $\mathbf{A}$ be a $\Sigma$-algebra. Then the \emph{support of} $\mathbf{A}$, denoted by $\mathrm{supp}_{S}(\mathbf{A})$, is $\mathrm{supp}_{S}(A)$, i.e., the support of the underlying $S$-sorted set $A$ of $\mathbf{A}$.
\end{definition}

\begin{remark}
The set $\{\mathrm{supp}_{S}(\mathbf{A})\mid \mathbf{A}\in \mathrm{Alg}(\Sigma)\}$ is a closure system on $S$.
\end{remark}

\begin{definition}
Let $\mathbf{A}$ be a $\Sigma$-algebra. We will say that $\mathbf{A}$ is \emph{finite} if $A$, the underlying $S$-sorted set of $\mathbf{A}$, is finite.
\end{definition}

\begin{remark}
In $\mathsf{Alg}(\Sigma)$, as it was the case with $\mathsf{Set}^{S}$, there is another notion of finiteness: A $\Sigma$-algebra  $\mathbf{A}$ is called $S$-\emph{finite} or \emph{locally finite}, abbreviated as $S\text{-}\mathrm{f}$, if and only if the underlying $S$-sorted set of $\mathbf{A}$ is $S$-finite. As was noted above this notion of finiteness plays a relevant role in the field of many-sorted algebra, e.g., to define $S$-finite, also called locally finite, terms and to distinguish, on the one hand, between many-sorted varieties and finitary many-sorted  varieties, and, on the other hand, between many-sorted quasivarieties and finitary many-sorted quasivarieties (see, e.g., \cite{gm85} and \cite{{m76}}).
\end{remark}

\section{Subalgebras}\label{SAlgSub}
We shall now go on to define the notion of subalgebra of a $\Sigma$-algebra $\mathbf{A}$, the principle of proof by Algebraic Induction, and the subalgebra generating operator for $\mathbf{A}$. The constructive description of the subalgebra generating operator for $\mathbf{A}$ will allow us to introduce the Borel subsets of a many-sorted algebra and the height of an element with respect to a given subset. This latter notion will be compared with the subelement preorder on $\mathbf{A}$.

\begin{definition}\label{DAlgSub}
Let $\mathbf{A}$ be a $\Sigma$-algebra and $X\subseteq A$. Given $(\mathbf{s},s)\in S^{\star}\times S$ and $\sigma\in\Sigma_{\mathbf{s},s}$, we will say that $X$ is \emph{closed under the operation} $\sigma^{\mathbf{A}}\colon A_{\mathbf{s}}\mor A_{s}$ if, for every $(a_{j})_{j\in\bb{\mathbf{s}}}\in X_{\mathbf{s}}$, $\sigma^{\mathbf{A}}((a_{j})_{j\in\bb{\mathbf{s}}})\in X_{s}$. We will say that $X$ is a \emph{closed subset} of $\mathbf{A}$ if $X$ is closed under the operations of $\mathbf{A}$. We will denote by $\mathrm{Cl}(\mathbf{A})$ the set of all closed subsets of $\mathbf{A}$ (which is an algebraic closure system on $A$) and by  $\mathbf{Cl}(\mathbf{A})$ the algebraic lattice $(\mathrm{Cl}(\mathbf{A}),\subseteq)$. We will say that a $\Sigma$-algebra $\mathbf{B}$ is a \emph{subalgebra} of $\mathbf{A}$ if $B\subseteq A$ and the canonical embedding of $B$ into $A$ determines an embedding of $\mathbf{B}$ into $\mathbf{A}$. We will denote by $\mathrm{Sub}(\mathbf{A})$ the set of all subalgebras of $\mathbf{A}$. Since $\mathrm{Cl}(\mathbf{A})$ and $\mathrm{Sub}(\mathbf{A})$ are isomorphic, we shall feel free to deal either with a closed subset of $\mathbf{A}$ or with the correlated subalgebra of $\mathbf{A}$, whichever is most convenient for the work at hand.
\end{definition}

\begin{definition}
Let $\mathbf{A}$ be a $\Sigma$-algebra. Then we denote by $\mathrm{Sg}_{\mathbf{A}}$ the many-sorted closure operator on $A$ defined as follows:
$$\textstyle
  \mathrm{Sg}_{\mathbf{A}}\nfunction
  {\mathrm{Sub}(A)}
  {\mathrm{Sub}(A)}
  {X}
  {\bigcap \{\,C\in\mathrm{Sub}(\mathbf{A})\mid X\subseteq C\,\}}
$$

We call $\mathrm{Sg}_{\mathbf{A}}$ the \emph{subalgebra generating many-sorted operator on} $A$ \emph{determined by} $\mathbf{A}$. For every $X\subseteq A$, we call $\mathrm{Sg}_{\mathbf{A}}(X)$ the \emph{subalgebra of} $\mathbf{A}$ \emph{generated by} $X$. Moreover, if $X\subseteq A$ is such that  $\mathrm{Sg}_{\mathbf{A}}(X) = A$, then we say that $X$ is an $S$-sorted set of \emph{generators} of $\mathbf{A}$, or that $X$ \emph{generates} $\mathbf{A}$. Besides, we say that $\mathbf{A}$ is \emph{finitely generated} if there exists an $S$-sorted  subset $X$ of $A$ such that $X$ \emph{generates} $\mathbf{A}$ and $\mathrm{card}(X)<\aleph_{0}$.
\end{definition}

\begin{remark}
Let $\mathbf{A}$ be a $\Sigma$-algebra. Then the algebraic closure operator $\mathrm{Sg}_{\mathbf{A}}$ is uniform, i.e., for every  $X$, $Y\subseteq A$, if $\mathrm{supp}_{S}(X) = \mathrm{supp}_{S}(Y)$, then we have $\mathrm{supp}_{S}(\mathrm{Sg}_{\mathbf{A}}(X)) = \mathrm{supp}_{S}(\mathrm{Sg}_{\mathbf{A}}(Y))$.
\end{remark}

We next recall the Principle of Proof by Algebraic Induction. 

\begin{proposition}[Principle of Proof by Algebraic Induction]\label{PPAI}
Let $\mathbf{A}$ be a  $\Sigma$-alge\-bra generated by $X$. Then to prove that a subset $Y$ of $A$ is equal to $A$ it suffices to show: (1) $X\subseteq Y$ (algebraic induction basis) and (2) $Y$ is a subalgebra of $\mathbf{A}$ (algebraic induction step).
\end{proposition}

\begin{proposition}\label{PAlgSg}
Let $\mathbf{A}$ be a $\Sigma$-algebra. Then the many-sorted closure operator $\mathrm{Sg}_{\mathbf{A}}$ on $A$ is algebraic, i.e., for every $S$-sorted subset $X$ of $A$, $\mathrm{Sg}_{\mathbf{A}}(X) =
\bigcup_{K\subseteq_{\mathrm{f}} X}\mathrm{Sg}_{\mathbf{A}}(K)$.
\end{proposition}

For a $\Sigma$-algebra $\mathbf{A}$ we next provide another, more constructive, description of the algebraic many-sorted closure operator $\Sg_{\mathbf{A}}$, which will allow us to define, for a $\Sigma$-algebra $\mathbf{A}$ and a subset $X$ of $A$, the Borel sets of $X$ in $\mathbf{A}$.

\begin{definition}
Let $\Sigma$ be an $S$-sorted signature and $\mathbf{A}$ a $\Sigma$-algebra.
\begin{enumerate}
\item We denote by $\mathrm{E}_{\mathbf{A}}$ the many-sorted operator on
      $A$ that assigns to an $S$-sorted subset $X$ of
      $A$, $\mathrm{E}_{\mathbf{A}}(X) =
      X\cup\bigl(\;\bigcup_{\sigma\in\Sigma_{\neq\lambda,s}}
      \sigma^{\mathbf{A}}[X_{\mathrm{ar}(\sigma)}]\bigr)_{s\in S}$, where, for $s\in
      S$, $\Sigma_{\neq\lambda,s}$ is the set of all many-sorted formal
      operations $\sigma$ such that $\mathrm{ar}(\sigma)\in S^{\star}-\{\lambda\}$,
      $\mathrm{car}(\sigma) = s$, and if $\mathrm{ar}(\sigma) = \mathbf{s}$, then
      $X_{\mathrm{ar}(\sigma)} = \prod_{i\in \bb{\mathbf{s}}}X_{s_{i}}$.

\item If $X\subseteq A$, then we define the family
      $(\mathrm{E}^{n}_{\mathbf{A}}(X))_{n\in \mathbb{N}}$ in $\mathrm{Sub}(A)$,
      recursively, as follows:
      \begin{align*}
      \mathrm{E}_{\mathbf{A}}^{0}(X) &=
      \textstyle
      X\cup (\bigcup_{\sigma\in \Sigma_{\lambda,s}}\sigma^{\mathbf{A}}[A_{\lambda}])_{s\in S}\text{,} \\
      \mathrm{E}_{\mathbf{A}}^{n+1}(X) &=
      \mathrm{E}_{\mathbf{A}}(\mathrm{E}_{\mathbf{A}}^{n}(X)) \text{, $n\geq 0$.}
      \end{align*}
\item We denote by $\mathrm{E}_{\mathbf{A}}^{\omega}$ the many-sorted operator on
      $A$ that assigns to an $S$-sorted subset $X$ of
      $A$, $\mathrm{E}_{\mathbf{A}}^{\omega}(X) = \bigcup_{n\in
      \mathbb{N}}\mathrm{E}_{\mathbf{A}}^{n}(X)$.
\end{enumerate}
\end{definition}

\begin{remark}
Since, for $s\in S$ and $\sigma\in\Sigma_{\lambda,s}$, $\sigma^{\mathbf{A}}$ is the value of $F^{\mathbf{A}}_{\sigma}\colon A_{\lambda}\mor A_{s}$ at the unique element in $A_{\lambda}$, it follows  that
$$
\textstyle
\bigl(\;\bigcup_{\sigma\in\Sigma_{\lambda,s}}\sigma^{\mathbf{A}}[A_{\lambda}]\bigr)_{s\in S} =
(\{\sigma^{\mathbf{A}}\mid \sigma\in \Sigma_{\lambda,s}\})_{s\in S}
$$
\end{remark}

\begin{proposition}\label{PSgE}
Let $\mathbf{A}$ be a $\Sigma$-algebra and $X\subseteq A$, then $\Sg_{\mathbf{A}}(X) = \mathrm{E}_{\mathbf{A}}^{\omega}(X)$.
\end{proposition}

\begin{proof}
Let us first prove that $\mathrm{Sg}_{\mathbf{A}}(X)\subseteq \mathrm{E}_{\mathbf{A}}^{\omega}(X)$. For this, because, by definition, $\mathrm{Sg}_{\mathbf{A}}(X)$ is the smallest subalgebra of $\mathbf{A}$ containing $X$, it suffices to verify  that $\mathrm{E}_{\mathbf{A}}^{\omega}(X)$ contains $X$ and that $\mathrm{E}_{\mathbf{A}}^{\omega}(X)$ is a subalgebra of $\mathbf{A}$. 

Since, by definition, $\mathrm{E}_{\mathbf{A}}^{0}(X) = X\cup (\bigcup_{\sigma\in \Sigma_{\lambda,s}}\sigma^{\mathbf{A}}[A_{\lambda}])_{s\in S}$, we have that
$X\subseteq \mathrm{E}_{\mathbf{A}}^{\omega}(X)$. Moreover, (1) if $s\in S$ and $\sigma\in\Sigma_{\lambda,s}$, then, since 
$$
\textstyle
\mathrm{E}_{\mathbf{A}}^{0}(X) 
= X\cup 
(\bigcup_{\sigma\in \Sigma_{\lambda,s}}\sigma^{\mathbf{A}}[A_{\lambda}])_{s\in S},$$
we have that $\sigma^{\mathbf{A}}\in \mathrm{E}_{\mathbf{A}}^{\omega}(X)_{s}$, and (2) if $\mathbf{s}\in S^{\star}-\{\lambda\}$, $s\in S$, $\sigma\in \Sigma_{\mathbf{s},s}$, and $(a_{j})_{j\in\bb{\mathbf{s}}}\in \mathrm{E}_{\mathbf{A}}^{\omega}(X)_{\mathbf{s}}$, then, for every $j\in \bb{\mathbf{s}}$, there exists an $n_{j}\in \mathbb{N}$ such that $a_{j}\in \mathrm{E}_{\mathbf{A}}^{n_{j}}(X)_{s_{i}}$. But the family
$(\mathrm{E}^{n}_{\mathbf{A}}(X))_{n\in \mathbb{N}}$ is an ascending chain, thus there exists a $n\in \bb{\mathbf{s}}$ such that, for every $j\in \bb{\mathbf{s}}$, $\mathrm{E}_{\mathbf{A}}^{n_{j}}(X)\subseteq
\mathrm{E}_{\mathbf{A}}^{n}(X)$. Therefore, for every $j\in \bb{\mathbf{s}}$, $a_{j}\in \mathrm{E}_{\mathbf{A}}^{n_{k}}(X)_{s_{j}}$, hence ${\sigma}^{\mathbf{A}}((a_{j})_{j\in\bb{\mathbf{s}}})\in
\mathrm{E}_{\mathbf{A}}^{n+1}(X)_{s}$. Whence ${\sigma}^{\mathbf{A}}((a_{j})_{j\in\bb{\mathbf{s}}})\in \mathrm{E}_{\mathbf{A}}^{\omega}(X)_{s}$. Consequently, $\mathrm{Sg}_{\mathbf{A}}(X)\subseteq \mathrm{E}_{\mathbf{A}}^{\omega}(X)$.

We next prove, by finite induction, that $\mathrm{E}_{\mathbf{A}}^{\omega}(X)\subseteq
\mathrm{Sg}_{\mathbf{A}}(X)$. 

\textsf{Base step of the induction.}

Since, by definition, $\mathrm{E}_{\mathbf{A}}^{0}(X) =  X\cup (\bigcup_{\sigma\in \Sigma_{\lambda,s}}\sigma^{\mathbf{A}}[A_{\lambda}])_{s\in S}$, 
we, obviously, have that $\mathrm{E}_{\mathbf{A}}^{0}(X)\subseteq\mathrm{Sg}_{\mathbf{A}}(X)$. 

\textsf{Inductive step of the induction.}

Let us suppose that, for $n\geq 0$, the requirement $\mathrm{E}_{\mathbf{A}}^{n}(X)\subseteq \mathrm{Sg}_{\mathbf{A}}(X)$ is satisfied. Then, since, by definition,   $\mathrm{E}_{\mathbf{A}}^{n+1}(X) = \mathrm{E}_{\mathbf{A}}(\mathrm{E}_{\mathbf{A}}^{n}(X))$, to prove that
$\mathrm{E}_{\mathbf{A}}^{n+1}(X)\subseteq \mathrm{Sg}_{\mathbf{A}}(X)$, it suffices to verify that $\mathrm{E}_{\mathbf{A}}^{n}(X)\subseteq\mathrm{Sg}_{\mathbf{A}}(X)$ and that 
$$
\textstyle
(\bigcup_{\sigma\in\Sigma_{\neq\lambda,s}}
\sigma^{\mathbf{A}}[\mathrm{E}_{\mathbf{A}}^{n}(X)_{\mathrm{ar}(\sigma)}])_{s\in S}
\subseteq\mathrm{Sg}_{\mathbf{A}}(X).$$

That $\mathrm{E}_{\mathbf{A}}^{n}(X)\subseteq\mathrm{Sg}_{\mathbf{A}}(X)$ follows directly from the induction hypothesis. Let $(\mathbf{s}, s)$ be an element of $(S^{\star}-\{\lambda\})\times S$, $\sigma$ an operation symbol in $\Sigma_{\mathbf{s},s}$, and $(a_{j})_{j\in\bb{\mathbf{s}}}\in
\mathrm{E}_{\mathbf{A}}^{n}(X)_{\mathbf{s}}$. Then,
$(a_{j})_{j\in\bb{\mathbf{s}}}\in \mathrm{Sg}_{\mathbf{A}}(X)_{\mathbf{s}}$, thus ${\sigma}^{\mathbf{A}}((a_{j})_{j\in\bb{\mathbf{s}}})\in \mathrm{Sg}_{\mathbf{A}}(X)_{s}$. Therefore $\sigma^{\mathbf{A}}[\mathrm{E}_{\mathbf{A}}^{n}(X)_{\mathrm{ar}(\sigma)}]\subseteq\mathrm{Sg}_{\mathbf{A}}(X)_{s}$. Hence $\bigcup_{\sigma\in\Sigma_{\neq\lambda,s}}\sigma^{\mathbf{A}}	[\mathrm{E}_{\mathbf{A}}^{n}(X)_{\mathrm{ar}(\sigma)}]\subseteq \mathrm{Sg}_{\mathbf{A}}(X)_{s}$. Consequently, $\mathrm{E}_{\mathbf{A}}^{\omega}(X)\subseteq\mathrm{Sg}_{\mathbf{A}}(X)$.
\end{proof}

The last proposition allows us to define the height of an element of type $s\in S$ in $\mathrm{Sg}_{\mathbf{A}}(X)_{s}$ as the smallest $n\in \mathbb{N}$ for which the element belongs to $\mathrm{E}_{\mathbf{A}}^{n}(X)_{s}$.

\begin{definition} Let $\mathbf{A}$ be a $\Sigma$-algebra, $X\subseteq A$, $s$ a sort in $S$, and $a$ an element in $\mathrm{Sg}_{\mathbf{A}}(X)_{s}$. Then we define the \emph{height of $a$ with respect to $\mathbf{A}$ and $X$}, denoted by $\mathrm{hgt}^{\mathbf{A}, X}_{s}(a)$, as:
$$
\mathrm{hgt}^{\mathbf{A}, X}_{s}(a) = \mathrm{min}\{n\in \mathbb{N}\mid a\in \mathrm{E}_{\mathbf{A}}^{n}(X)_{s}\}.
$$
\end{definition}

\begin{remark}
Let $\mathbf{A}$ be a $\Sigma$-algebra and $X\subseteq A$. Then $\mathrm{hgt}^{\mathbf{A},X} = (\mathrm{hgt}^{\mathbf{A},X}_{s})_{s\in S}$ is an $S$-sorted mapping from $\mathrm{Sg}_{\mathbf{A}}(X)$ to $(\mathbb{N})_{s\in S}$.
\end{remark}


We next define, for a $\Sigma$-algebra $\mathbf{A}$ and a subset $X\subseteq A$, the family of the Borel sets of $X$ in $\mathbf{A}$.

\begin{definition} Let $\mathbf{A}$ be a $\Sigma$-algebra and $X\subseteq A$. Then we define the family $(\mathrm{B}^{n}_{\mathbf{A}}(X))_{n\in \mathbb{N}}$ in $\mathrm{Sub}(A)$, recursively, as follows:
      \begin{align*}
      \mathrm{B}_{\mathbf{A}}^{0}(X) &= \mathrm{E}^{0}_{\mathbf{A}}(X)\text{,} \\
       \mathrm{B}_{\mathbf{A}}^{n+1}(X) &=
      \mathrm{E}_{\mathbf{A}}^{n+1}(X) -  \mathrm{E}_{\mathbf{A}}^{n}(X)\text{, $n\geq 0$.}
      \end{align*}
For a natural number $n\in\mathbb{N}$, the $S$-sorted set $\mathrm{B}^{n}_{\mathbf{A}}(X)$ will be called the \emph{$n$-th strict Borel set of $X$ in $\mathbf{A}$}, abbreviated to \emph{$n$-th Borel set of $X$ in $\mathbf{A}$}. 
\end{definition}

\begin{remark}
Let $\mathbf{A}$ be a $\Sigma$-algebra and $X\subseteq A$. Then  
$$
\textstyle
\mathrm{B}_{\mathbf{A}}^{0}(X) = X\cup (\bigcup_{\sigma\in \Sigma_{\lambda,s}}\sigma^{\mathbf{A}}[A_{\lambda}])_{s\in S},
$$ 
while the $1$-th Borel set of $X$ in $\mathbf{A}$ is
$$
\textstyle
\mathrm{B}_{\mathbf{A}}^{1}(X) = \bigl(\;\bigcup_{\tau\in\Sigma_{\neq\lambda,s}}
      \tau^{\mathbf{A}}[(X\cup (\bigcup_{\sigma\in \Sigma_{\lambda,s}}\sigma^{\mathbf{A}}[A_{\lambda}])_{s\in S})_{\mathrm{ar}(\sigma)}]\bigr)_{s\in S}.
$$
\end{remark}

We now characterize, for every natural number $n$, the $n$-th Borel set of a set $X$ in a $\Sigma$-algebra $\mathbf{A}$ by means of the $S$-sorted mapping $\mathrm{hgt}^{\mathbf{A},X}$.

\begin{proposition} Let $\mathbf{A}$ be a $\Sigma$-algebra, $X\subseteq A$, $n\in \mathbb{N}$ and $s\in S$. Then it happens that:
$$
\mathrm{B}^{n}_{\mathbf{A}}(X)_{s}
=
\{a\in\mathrm{Sg}_{\mathbf{A}}(X)_{s}\mid \mathrm{hgt}^{\mathbf{A},X}_{s}(a)=n\}.
$$
Therefore, for every $n\in \mathbb{N}$, $\mathrm{B}^{n}_{\mathbf{A}}(X) = ((\mathrm{hgt}^{\mathbf{A},X}_{s})^{-1}[\{n\}])_{s\in S}$.
\end{proposition}

We next define for a $\Sigma$-algebra $\mathbf{A}$ several binary relations on 
$\coprod A$.

\begin{definition}\label{DAlgOrd} 
Let $\mathbf{A}$ be a $\Sigma$-algebra.
We let $\prec_{\mathbf{A}}$ denote the binary relation on $\coprod A$ consisting of the ordered pairs $((b,t),(a,s))\in(\coprod A)^{2}$ for which there exists a word $\mathbf{s}\in S^{\star}-\{\lambda\}$, an operation symbol $\sigma\in \Sigma_{\mathbf{s},s}$ and a family of elements $(a_{j})_{j\in\bb{\mathbf{s}}}\in A_{\mathbf{s}}$ such that $a=\sigma^{\mathbf{A}}((a_{j})_{j\in\bb{\mathbf{s}}})$ and, for some $j\in\bb{\mathbf{s}}$, $s_{j}=t$ and $b=a_{j}$. If $((b,t),(a,s))\in \prec_{\mathbf{A}}$, also written as $(b,t)\prec_{\mathbf{A}}(a,s)$, then we will say that $(b,t)$ $\prec_{\mathbf{A}}$-\emph{precedes} $(a,s)$ or, when no confusion can arise, that $b$ \emph{precedes} $a$. We will call $\prec_{\mathbf{A}}$ the \emph{algebraic predecessor relation} on $\mathbf{A}$.

We will denote by $\leq_{\mathbf{A}}$ the reflexive and transitive closure of $\prec_{\mathbf{A}}$, i.e., the preorder on $\coprod A$ generated by $\prec_{\mathbf{A}}$, and by $<_{\mathbf{A}}$ the transitive closure of $\prec_{\mathbf{A}}$. If $((b,t),(a,s))\in \leq_{\mathbf{A}}$, also written as $(b,t)\leq_{\mathbf{A}}(a,s)$, then we will say that $b$ is a \emph{subelement} of $a$. The preorder $\leq_{\mathbf{A}}$ will be called the \emph{algebraic predecessor preorder} on $\mathbf{A}$.
\end{definition}

\begin{remark}\label{RAlgOrd}
The preorder $\leq_{\mathbf{A}}$ on $\coprod A$ is $\bigcup_{n\in\mathbb{N}}\prec_{\mathbf{A}}^{n}$, where $\prec_{\mathbf{A}}^{0}$ is the diagonal of $\coprod A$ and, for $n\in\mathbb{N}$, 
$$\prec_{\mathbf{A}}^{n+1} = \prec_{\mathbf{A}}^{n}\circ \prec_{\mathbf{A}}.$$ 
Thus, for every $((b,t),(a,s))\in (\coprod A)^{2}$, $(b,t)\leq_{\mathbf{A}}(a,s)$ if and only if $s = t$ and $b = a$ or there exists  a natural number $m\in\mathbb{N}-\{0\}$, a word $\mathbf{c}\in S^{\star}$ of length $\bb{\mathbf{c}}=m+1$, and a family of elements $(r_{k})_{k\in\bb{\mathbf{c}}}$ in $\mathrm{A}_{\mathbf{c}}$, such that $c_{0}=t$,  $r_{0}=b$, $c_{m}=s$, $r_{m}=a$ and, for every 
$k\in m$, $(r_{k},c_{k})\prec_{\mathrm{A}}(r_{k+1},c_{k+1})$. Moreover, $<_{\mathbf{A}}$, the transitive closure of $\prec_{\mathbf{A}}$, is $\bigcup_{n\in\mathbb{N}-
\{0\}}\prec_{\mathbf{A}}^{n}$.
\end{remark}

We conclude this subsection by showing how $\prec_{\mathbf{A}}$ relates to the height of an element of the underlying $S$-sorted set of $\mathbf{A}$.

\begin{proposition} 
Let $\mathbf{A}$ be a $\Sigma$-algebra, $s,t$ sorts in $S$ and $a,b$ elements in $A_{s}$ and $A_{t}$, respectively, such that $(b,t)$ $\prec_{\mathbf{A}}$-precedes $(a,s)$, thus there exists a word 
$\mathbf{s}\in S^{\star}-\{\lambda\}$, an operation symbol $\sigma\in \Sigma_{\mathbf{s},s}$ and a family of elements $(a_{j})_{j\in\bb{\mathbf{s}}}\in A_{\mathbf{s}}$ such that $a=\sigma^{\mathbf{A}}((a_{j})_{j\in\bb{\mathbf{s}}})$ and, for some $j\in\bb{\mathbf{s}}$, $s_{j}=t$ and $b=a_{j}$. Then, for every $S$-sorted subset $X\subseteq A$ it holds that
$$
\mathrm{hgt}^{\mathbf{A},X}_{s}(a)\geq
\mathrm{max}\{
\mathrm{hgt}^{\mathbf{A},X}_{s_{j}}(a_{j})\mid j\in\bb{\mathbf{s}}
\}+1.
$$
In particular, 
$
\mathrm{hgt}^{\mathbf{A},X}_{s}(a)\geq
\mathrm{hgt}^{\mathbf{A},X}_{t}(b)+1
$.
\end{proposition}

\section{Congruences}\label{SCong}
Our next goal is to define the concepts of congruence on a $\Sigma$-algebra and of quotient of a $\Sigma$-algebra by a congruence on it. Moreover, we recall the notion of kernel of a homomorphism between $\Sigma$-algebras, the universal property of the quotient of a $\Sigma$-algebra by a congruence on it, and, for a $\Sigma$-algebra, the congruence generating operator for it. 

\begin{definition}\label{DCong}
Let $\mathbf{A}$ be a $\Sigma$-algebra and $\Phi$ an $S$-sorted equivalence on $A$. We will say that $\Phi$ is an $S$-\emph{sorted congruence on} $\mathbf{A}$ if, for every $(\mathbf{s},s)\in (S^{\star}-\{\lambda\})\times S$, every $\sigma\in \Sigma_{\mathbf{s},s}$, and every $(a_{j})_{j\in\bb{\mathbf{s}}},(b_{j})_{j\in\bb{\mathbf{s}}}\in A_{\mathbf{s}}$, if, for every $j\in \bb{\mathbf{s}}$, $(a_{j}, b_{j})\in\Phi_{s_{j}}$, then $(\sigma^{\mathbf{A}}((a_{j})_{j\in\bb{\mathbf{s}}}), \sigma^{\mathbf{A}}((b_{j})_{j\in\bb{\mathbf{s}}}))\in \Phi_{s}$.
Sometimes, we will abbreviate the expression ``$S$-sorted congruence on $\mathbf{A}$'' to ``congruence on $\mathbf{A}$''. We will denote by $\mathrm{Cgr}(\mathbf{A})$ the set of all $S$-sorted congruences on $\mathbf{A}$ (which is an algebraic closure system on $A\times A$), by $\mathbf{Cgr}(\mathbf{A})$ the algebraic lattice  $(\mathrm{Cgr}(\mathbf{A}),\subseteq)$, by $\nabla_{\mathbf{A}}$ the greatest element of $\mathbf{Cgr}(\mathbf{A})$, and by $\Delta_{\mathbf{A}}$ the least element of $\mathbf{Cgr}(\mathbf{A})$.

For a congruence $\Phi$ on $\mathbf{A}$, the \emph{quotient $\Sigma$-algebra of} $\mathbf{A}$ \emph{by} $\Phi$, denoted by $\mathbf{A}/\Phi$, is the $\Sigma$-algebra $(A/\Phi, F^{\mathbf{A}/\Phi})$, where, for every  $(\mathbf{s},s)\in S^{\star}\times S$ and every $\sigma\in \Sigma_{\mathbf{s},s}$, the operation $F^{\mathbf{A}/\Phi}_{\sigma}$ from $(A/\Phi)_{\mathbf{s}}$ to $A_{s}/\Phi_{s}$, also denoted, to simplify, by $\sigma^{\mathbf{A}/\Phi}$, sends $([a_{j}]_{\Phi_{s_{j}}})_{j\in\bb{\mathbf{s}}}$ in $(A/\Phi)_{\mathbf{s}}$ to $[\sigma^{\mathbf{A}}((a_{j})_{j\in \bb{\mathbf{s}}})]_{\Phi_{s}}$ in $A_{s}/\Phi_{s}$,  %
and the \emph{canonical projection} from $\mathbf{A}$ to $\mathbf{A}/\Phi$, denoted by $\mathrm{pr}^{\Phi}\colon \mathbf{A}\mor \mathbf{A}/\Phi$, is the full surjective homomorphism determined by the projection from $A$ to $A/\Phi$.
The ordered pair $(\mathbf{A}/\Phi,\mathrm{pr}^{\Phi})$ has the following universal property: $\mathrm{Ker}(\mathrm{pr}^{\Phi})$ is $\Phi$ and, for every $\Sigma$-algebra $\mathbf{B}$ and every homomorphism $f$ from $\mathbf{A}$ to $\mathbf{B}$, if $\Phi\subseteq \mathrm{Ker}(f)$, then there exists a unique homomorphism $h$ from $\mathbf{A}/\Phi$ to $\mathbf{B}$ such that $h\circ\mathrm{pr}^{\Phi} = f$. In particular, if $\Psi$ is a congruence on $A$ such that $\Phi\subseteq \Psi$, then we will denote by $\mathrm{p}^{\Phi,\Psi}$ the unique homomorphism from $\mathbf{A}/\Phi$ to $\mathbf{A}/\Psi$ such that $\mathrm{p}^{\Phi,\Psi}\circ \mathrm{pr}^{\Phi} = \mathrm{pr}^{\Psi}$. 
\end{definition}


\begin{remark}
Let $\mathsf{ClfdAlg}(\Sigma)$ be the category whose objects are the \emph{classified $\Sigma$-algebras}, i.e, the ordered pairs $(\mathbf{A},\Phi)$ where $\mathbf{A}$ is a $\Sigma$-algebra and $\Phi$ a congruence on $\mathbf{A}$, and in which the set of morphisms from $(\mathbf{A},\Phi)$ to $(\mathbf{B},\Psi)$ is the set of all homomorphisms $f$ from $\mathbf{A}$ to $\mathbf{B}$ such that, for every $s\in S$ and every $(x,y)\in A^{2}_{s}$, if $(x,y)\in \Phi_{s}$, then $(f_{s}(x),f_{s}(y))\in \Psi_{s}$. Let $G$ be the functor from $\mathsf{Alg}(\Sigma)$ to $\mathsf{ClfdAlg}(\Sigma)$ whose object mapping sends $\mathbf{A}$ to $(\mathbf{A},\Delta_{\mathbf{A}})$ and whose morphism mapping sends $f\colon \mathbf{A}\mor \mathbf{B}$ to $f\colon (\mathbf{A},\Delta_{\mathbf{A}})\mor (\mathbf{B},\Delta_{\mathbf{B}})$. Then, for every classified $\Sigma$-algebra $(\mathbf{A},\Phi)$, there exists a universal mapping from $(\mathbf{A},\Phi)$ to $G$, which is precisely the ordered pair $(\mathbf{A}/\Phi,\mathrm{pr}^{\Phi})$ with $\mathrm{pr}^{\Phi}\colon (\mathbf{A},\Phi)\mor (\mathbf{A}/\Phi,\Delta_{\mathbf{A}/\Phi})$.
\end{remark}

We next recall the congruence generating operator for a $\Sigma$-algebra and provide a constructive description for it.

\begin{definition}\label{DCongOpInt} Let $\mathbf{A}$ be a $\Sigma$-algebra. Then we denote by $\mathrm{Cg}_{\mathbf{A}}$ the many-sorted closure operator on $A\times A$ defined as follows:

$$\textstyle
  \mathrm{Cg}_{\mathbf{A}}\nfunction
  {\mathrm{Sub}(A\times A)}
  {\mathrm{Sub}(A\times A)}
  {\Phi}
  {\bigcap \{\,\Psi\in\mathrm{Con}(\mathbf{A})\mid \Phi\subseteq \Psi\,\}}
$$

We call $\mathrm{Cg}_{\mathbf{A}}$ the \emph{congruence generating many-sorted operator on} $A\times A$ \emph{determined by} $\mathbf{A}$. For every relation $\Phi\subseteq A\times A$, we call $\mathrm{Cg}_{\mathbf{A}}(\Phi)$ the \emph{congruence on} $\mathbf{A}$ \emph{generated by} $\Phi$. 
\end{definition}

For a $\Sigma$-algebra $\mathbf{A}$ we next provide another, more constructive, description of the algebraic many-sorted closure operator $\mathrm{Cg}_{\mathbf{A}}$.

\begin{definition}\label{DCongOpC} 
Let $\Sigma$ be an $S$-sorted signature and $\mathbf{A}$ a $\Sigma$-algebra.
\begin{itemize}
\item[(1)] We denote by $\mathrm{C}_{\mathbf{A}}$ the many-sorted operator on $A\times A$ that assigns to an $S$-sorted relation $\Phi\subseteq A\times A$, the $S$-sorted relation 
$$
\textstyle
\mathrm{C}_{\mathbf{A}}(\Phi)
=(\Phi\circ\Phi)
\cup
\left(
\bigcup_{\sigma\in\Sigma_{\neq\lambda,s}}
\sigma^{\mathbf{A}}\times \sigma^{\mathbf{A}}
\left[
\Phi_{\mathrm{ar}(\sigma)}
\right]
\right)_{s\in S},
$$
where, for $s\in S$, $\Sigma_{\neq\lambda,s}$ is the set of all many-sorted formal operations $\sigma$ such that $\mathrm{ar}(\sigma)\in S^{\star}-\{\lambda\}$, $\mathrm{car}(\sigma)=s$, and if $\mathrm{ar}(\sigma)=\mathbf{s}$, then $\Phi_{\mathrm{ar}(\sigma)}=\prod_{j\in\bb{\mathbf{s}}}\Phi_{s_{j}}$. Let us note that 
$$
\textstyle
\Phi_{\mathrm{ar}(\sigma)}\subseteq \prod_{j\in\bb{\mathbf{s}}}A_{s_{j}}^{2}\cong (\prod_{j\in\bb{\mathbf{s}}}A_{s_{j}})^{2}.
$$
\item[(2)] If $\Phi\subseteq A\times A$, then we define the family $(\mathrm{C}^{n}_{\mathbf{A}}(\Phi))_{n\in\mathbb{N}}$ in $\mathrm{Sub}(A\times A)$, recursively, as follows:
\begin{align*}
\mathrm{C}^{0}_{\mathbf{A}}(\Phi)&=\Phi\cup\Phi^{-1}\cup\Delta_{\mathrm{A}},\\
\mathrm{C}^{n+1}_{\mathbf{A}}(\Phi)&=\mathrm{C}_{\mathbf{A}}(\mathrm{C}_{\mathbf{A}}^{n}(\Phi)),\, n\geq 0.
\end{align*}
\item[(3)] We denote by $\mathrm{C}^{\omega}_{\mathbf{A}}$ the many-sorted operator on $A\times A$ that assigns to an $S$-sorted relation $\Phi\subseteq A\times A$, $\mathrm{C}^{\omega}_{\mathbf{A}}(\Phi)=\bigcup_{n\in\mathbb{N}}\mathrm{C}^{n}_{\mathbf{A}}(\Phi)$.
\end{itemize}
\end{definition}

\begin{proposition}\label{PCongOpC} 
Let $\mathbf{A}$ be a $\Sigma$-algebra and $\Phi\subseteq A\times A$ a relation on $A$. Then, for every $n\in\mathbb{N}$, we have that 
\begin{itemize}
\item[(1)] $\mathrm{C}^{n}_{\mathbf{A}}(\Phi)$ is a reflexive relation;
\item[(2)]  $\mathrm{C}^{n}_{\mathbf{A}}(\Phi)\subseteq\mathrm{C}^{n+1}_{\mathbf{A}}(\Phi)$;
\item[(3)] $\mathrm{C}^{n}_{\mathbf{A}}(\Phi)$ is a symmetric relation.
\end{itemize}
\end{proposition}

\begin{proof}
We prove the first item by induction on $n\in\mathbb{N}$.

\textsf{Base step of the induction for the reflexivity property.}

Since $\Delta_{\mathrm{A}}\subseteq \mathrm{C}^{0}_{\mathbf{A}}(\Phi)$, we have that $\mathrm{C}^{0}_{\mathbf{A}}(\Phi)$ is a reflexive relation.

\textsf{Inductive step of the induction  for the reflexivity property.}

Let us assume that the statement holds for $n\in\mathbb{N}$, i.e., that $\mathrm{C}^{n}_{\mathbf{A}}(\Phi)$ is a reflexive relation.

We now prove that the statement holds for $n+1$. By the inductive hypothesis, $\mathrm{C}^{n}_{\mathbf{A}}(\Phi)$ is a reflexive relation, i.e., $\Delta_{\mathrm{A}}\subseteq \mathrm{C}^{n}_{\mathbf{A}}(\Phi)$. Hence 
$$
\Delta_{\mathrm{A}}\circ \mathrm{C}^{n}_{\mathbf{A}}(\Phi)\subseteq \mathrm{C}^{n}_{\mathbf{A}}(\Phi)\circ \mathrm{C}^{n}_{\mathbf{A}}(\Phi) \quad\text{and}\quad \Delta_{\mathrm{A}}\circ \Delta_{\mathrm{A}}\subseteq \Delta_{\mathrm{A}}\circ \mathrm{C}^{n}_{\mathbf{A}}(\Phi).
$$ 
But $\Delta_{\mathrm{A}}\subseteq \Delta_{\mathrm{A}}\circ \Delta_{\mathrm{A}}$. Consequently $\Delta_{\mathrm{A}}\subseteq \mathrm{C}^{n}_{\mathbf{A}}(\Phi)\circ \mathrm{C}^{n}_{\mathbf{A}}(\Phi)$. Moreover, we have that $\mathrm{C}^{n}_{\mathbf{A}}(\Phi)\circ \mathrm{C}^{n}_{\mathbf{A}}(\Phi)\subseteq \mathrm{C}^{n+1}_{\mathbf{A}}(\Phi)$. Thus $\Delta_{\mathrm{A}}\subseteq\mathrm{C}^{n+1}_{\mathbf{A}}(\Phi)$

Therefore, $\mathrm{C}^{n+1}_{\mathbf{A}}(\Phi)$ is a reflexive relation. This finishes the proof of the first item.

We now prove the second item. By (1), for every $n\in\mathbb{N}$, the relation $\mathrm{C}^{n}_{\mathbf{A}}(\Phi)$ is reflexive. Hence $\Delta_{\mathrm{A}}\subseteq \mathrm{C}^{n}_{\mathbf{A}}(\Phi)$. Since $\Delta_{\mathrm{A}}$ is a neutral element for the composition of relations, we have that 
$$
\mathrm{C}^{n}_{\mathbf{A}}(\Phi)
=
\mathrm{C}^{n}_{\mathbf{A}}(\Phi)\circ\Delta_{\mathrm{A}}
\subseteq \mathrm{C}^{n}_{\mathbf{A}}(\Phi)\circ\mathrm{C}_{\mathbf{A}}^{n}(\Phi)
\subseteq \mathrm{C}^{n+1}_{\mathbf{A}}(\Phi).
$$
Hence, $\mathrm{C}^{n}_{\mathbf{A}}(\Phi)\subseteq\mathrm{C}^{n+1}_{\mathbf{A}}(\Phi)$. This finishes the proof of the second item.

We now prove the third statement by induction on $\mathrm{n}\in\mathbb{N}$.

\textsf{Base step of the induction for the symmetry property.}

Let us recall that 
$\mathrm{C}^{0}_{\mathbf{A}}(\Phi)=\Phi\cup\Phi^{-1}\cup\Delta_{\mathrm{A}}$. Hence, the inverse of  $\mathrm{C}^{0}_{\mathbf{A}}(\Phi)$, which is given by 
$
\left(
\mathrm{C}^{0}_{\mathbf{A}}(\Phi)
\right)^{-1}
=\Phi^{-1}\cup\Phi\cup\Delta_{\mathrm{A}}
$, is included in $\mathrm{C}^{0}_{\mathbf{A}}(\Phi)$. Hence, $\mathrm{C}^{0}_{\mathbf{A}}(\Phi)$ is a symmetric relation.

\textsf{Inductive step of the induction for the symmetry property.}

Let us assume that the statement holds for $n\in\mathbb{N}$, i.e., that $\mathrm{C}^{n}_{\mathbf{A}}(\Phi)$ is a symmetric relation.

We now prove the statement for $n+1$. Let $s$ be a sort in $s$ and let  $(x,y)$ be a pair in $\mathrm{C}^{n+1}_{\mathbf{A}}(\Phi)$. Let us recall that $\mathrm{C}^{n+1}_{\mathbf{A}}(\Phi)$ is given by
$$
\textstyle
\mathrm{C}^{n+1}_{\mathbf{A}}(\Phi)=
\mathrm{C}_{\mathbf{A}}(\mathrm{C}^{n}_{\mathbf{A}}(\Phi))
=
(\mathrm{C}^{n}_{\mathbf{A}}(\Phi)\circ\mathrm{C}^{n}_{\mathbf{A}}(\Phi))
\cup
\left(
\bigcup_{\sigma\in\Sigma_{\neq\lambda,s}}
\sigma^{\mathbf{A}}\times \sigma^{\mathbf{A}}
\left[
\mathrm{C}^{n}_{\mathbf{A}}(\Phi)
\right]
\right)_{s\in S},
$$
hence the pair $(x,y)$ is either (1) a pair  in $\mathrm{C}^{n}_{\mathbf{A}}(\Phi)_{s}\circ\mathrm{C}^{n}_{\mathbf{A}}(\Phi)_{s}$ or (2) a pair in $\bigcup_{\sigma\in\Sigma_{\neq\lambda,s}}
\sigma^{\mathbf{A}}\times \sigma^{\mathbf{A}}
\left[
\mathrm{C}^{n}_{\mathbf{A}}(\Phi)
\right]
$.

If (1), then exists an element $z\in A_{s}$ such that the pairs $(x,z)$ and $(z,y)$ are in $\mathrm{C}^{n}_{\mathbf{A}}(\Phi)_{s}$. By induction, $\mathrm{C}^{n}_{\mathbf{A}}(\Phi)_{s}$ is a symmetric relation. Hence, the pairs $(z,x)$ and $(y,z)$ belong to $\mathrm{C}^{n}_{\mathbf{A}}(\Phi)_{s}$. Thus, the pair $(y,x)$ also belongs to $\mathrm{C}^{n+1}_{\mathbf{A}}(\Phi)_{s}$.

If (2), then we can find a non-empty word $\mathbf{s}\in S^{\star}-\{\lambda\}$, an operation symbol $\sigma\in \Sigma_{\mathbf{s},s}$ and a family of pairs $((x_{j},y_{j}))_{j\in\bb{\mathbf{s}}}$ in $\mathrm{C}^{n}_{\mathbf{A}}(\Phi)_{\mathbf{s}}$ such that $x=\sigma^{\mathbf{A}}((x_{j})_{j\in\bb{\mathbf{s}}})$ and $y=\sigma^{\mathbf{A}}((y_{j})_{j\in\bb{\mathbf{s}}})$. By induction we have that, for every $j\in\bb{\mathbf{s}}$, $\mathrm{C}^{n}_{\mathbf{A}}(\Phi)_{s_{j}}$ is a symmetric relation. Hence, for every $j\in\bb{\mathbf{s}}$, the pair $(y_{j},x_{j})$ belongs to $\mathrm{C}^{n}_{\mathbf{A}}(\Phi)_{s_{j}}$. Thus, the pair $(\sigma^{\mathbf{A}}((y_{j})_{j\in\bb{\mathbf{s}}}), \sigma^{\mathbf{A}}((x_{j})_{j\in\bb{\mathbf{s}}}))$, which is equal to $(y,x)$, also belongs to $\mathrm{C}^{n+1}_{\mathbf{A}}(\Phi)_{s}$.

In any case, the pair $(y,x)$ belongs to $\mathrm{C}^{n+1}_{\mathbf{A}}(\Phi)_{s}$.
Hence, $\mathrm{C}^{n+1}_{\mathbf{A}}(\Phi)$ is a symmetric relation. This finishes the proof of the third item.
\end{proof}

\begin{proposition}\label{PCongIntC} 
Let $\mathbf{A}$ be a $\Sigma$-algebra and $\Phi\subseteq A\times A$. Then
$\mathrm{Cg}_{\mathbf{A}}(\Phi)=\mathrm{C}_{\mathbf{A}}^{\omega}(\Phi)$.
\end{proposition}

\begin{proof}
We first prove the inclusion from left to right. To do so, we prove that $\mathrm{C}_{\mathbf{A}}^{\omega}(\Phi)$ is a $\Sigma$-congruence on $\mathbf{A}$.

By Proposition~\ref{PCongOpC}, $\Delta_{A}\subseteq\mathrm{C}^{0}_{\mathbf{A}}(\Phi)$, it follows that $\Delta_{A}\subseteq\bigcup_{n\in\mathbb{N}}\mathrm{C}^{n}_{\mathbf{A}}(\Phi)=\mathrm{C}^{\omega}_{\mathbf{A}}(\Phi)$. Hence, $\mathrm{C}^{\omega}_{\mathbf{A}}(\Phi)$ is a reflexive relation. 

Regarding symmetry, let $s$ be a sort in $S$ and let $(x,y)$ be a pair in $\mathrm{C}^{\omega}_{\mathbf{A}}(\Phi)_{s}$, then there exists an $n\in\mathbb{N}$ for which $(x,y)$ is a pair in $\mathrm{C}^{n}_{\mathbf{A}}(\Phi)_{s}$. By Proposition~\ref{PCongOpC}, $\mathrm{C}^{n}_{\mathbf{A}}(\Phi)_{s}$ is a symmetric relation. Hence $(y,x)\in \mathrm{C}^{n}_{\mathbf{A}}(\Phi)_{s}\subseteq \mathrm{C}^{\omega}_{\mathbf{A}}(\Phi)_{s}$. Therefore $\mathrm{C}^{\omega}_{\mathbf{A}}(\Phi)$ is a symmetric relation. 

Regarding transitivity, let $s$ be a sort in $S$ and let $(x,y)$ and $(y,z)$ be a pair of terms in $\mathrm{C}^{\omega}_{\mathbf{A}}(\Phi)_{s}$, then there exists an $m$ and an $n\in\mathbb{N}$ for which $(x,y)$ is a pair in $\mathrm{C}^{m}_{\mathbf{A}}(\Phi)_{s}$ and $(y,z)$ is a pair in $\mathrm{C}^{n}_{\mathbf{A}}(\Phi)_{s}$. Let $r=\mathrm{max}(m,n)$ then, in virtue of Proposition~\ref{PCongOpC}, the pairs $(x,y)$ and $(y,z)$ are in $\mathrm{C}^{r}_{\mathbf{A}}(\Phi)_{s}$. It follows that $(x,z)$ is a pair in $\mathrm{C}^{r}_{\mathbf{A}}(\Phi)_{s}\circ\mathrm{C}^{r}_{\mathbf{A}}(\Phi)_{s}\subseteq \mathrm{C}^{r+1}_{\mathbf{A}}(\Phi)_{s}\subseteq \mathrm{C}^{\omega}_{\mathbf{A}}(\Phi)_{s}$. Therefore $\mathrm{C}^{\omega}_{\mathbf{A}}(\Phi)$ is a transitive relation. 

Finally, we check that $\mathrm{C}^{m}_{\mathbf{A}}(\Phi)$ is compatible with the operations in $\Sigma$. Let $\mathbf{s}$ be a non-empty word in $S^{\star}-\{\lambda\}$ and $s$ a sort in $S$. Let $\sigma$ be an operation symbol in $\Sigma_{\mathbf{s},s}$. Let $(x_{j})_{j\in\bb{\mathbf{s}}}$ and $(y_{j})_{j\in\mathbf{s}}$ be a pair of families in $\mathbf{A}_{\mathbf{s}}$ such that, for every $j\in\bb{\mathbf{s}}$, the pair $(x_{j},y_{j})$ is in $\mathrm{C}^{\omega}_{\mathbf{A}}(\Phi)_{s_{j}}$. Then there exists, for every $j\in\bb{\mathbf{s}}$, a natural number $n_{j}\in\mathbb{N}$ such that $(x_{j},y_{j})$ is a pair in $\mathrm{C}^{n_{j}}_{\mathbf{A}}(\Phi)_{s_{j}}$. Let $r=\mathrm{max}\{n_{j}\mid j\in\bb{\mathbf{s}}\}$ then, in virtue of Proposition~\ref{PCongOpC}, for every $j\in\bb{\mathbf{s}}$, the pair $(x_{j},y_{j})$ is in $\mathrm{C}^{r}_{\mathbf{A}}(\Phi)_{s_{j}}$. Therefore  
$
(
\sigma^{\mathbf{A}}
(
(
x_{j}
)_{j\in\bb{\mathbf{s}}}
),
\sigma^{\mathbf{A}}
(
(
y_{j}
)_{j\in\bb{\mathbf{s}}}
)
)
$ is a pair in $\sigma^{\mathbf{A}}\times \sigma^{\mathbf{A}}[\mathrm{C}^{r}_{\mathbf{A}}(\Phi)_{\mathbf{s}}]$,
which is a subset of $\mathrm{C}^{r+1}_{\mathbf{A}}(\Phi)_{s}\subseteq \mathrm{C}^{\omega}_{\mathbf{A}}(\Phi)_{s}$. 

All in all, we can affirm that $\mathrm{C}^{\omega}_{\mathbf{A}}(\Phi)$ is a congruence on $\mathbf{A}$.

Moreover $\mathrm{C}^{\omega}_{\mathbf{A}}(\Phi)$ contains $\Phi$. This is so because $\Phi\subseteq \mathrm{C}^{0}_{\mathbf{A}}(\Phi)$. From these facts, we can affirm that $\mathrm{Cg}_{\mathbf{A}}(\Phi)\subseteq\mathrm{C}^{\omega}_{\mathbf{A}}(\Phi)$.

To prove the other inclusion, we proceed by induction on $n\in \mathbb{N}$.

\textsf{Base step of the induction.} 

Since $\mathrm{Cg}_{\mathbf{A}}(\Phi)$ is a congruence on $\mathbf{A}$ containing $\Phi$, we can affirm that $\Delta_{A}$, $\Phi$ and $\Phi^{-1}$ are subsets of $\mathrm{Cg}_{\mathbf{A}}(\Phi)$. Hence $\mathrm{C}^{0}_{\mathbf{A}}(\Phi)\subseteq \mathrm{Cg}_{\mathbf{A}}(\Phi)$.

\textsf{Inductive step of the induction.}

Let us assume that the statement holds for  $n\in\mathbb{N}$, i.e., we have that
$\mathrm{C}^{n}_{\mathbf{A}}(\Phi)\subseteq \mathrm{Cg}_{\mathbf{A}}(\Phi)$. 

We now prove the statement for $n+1$. 
By the inductive hypothesis, the relation $\mathrm{C}^{n}_{\Phi}$ is included in $\mathrm{Cg}_{\mathbf{A}}(\Phi)$. Since $\mathrm{Cg}_{\mathbf{A}}(\Phi)$ is a transitive relation, we can affirm that $(\mathrm{C}^{n}_{\Phi}\circ \mathrm{C}^{n}_{\Phi})$ is included in $\mathrm{Cg}_{\mathbf{A}}(\Phi)$. 

On the other hand, for every non-empty word $\mathbf{s}\in S-\{\lambda\}$, every sort $s\in S$, every operation symbol $\sigma$ in $\Sigma_{\mathbf{s},s}$, if $(x_{j},y_{j})_{j\in\bb{\mathbf{s}}}$ is a family of pairs in $\mathrm{C}^{n}_{\mathbf{A}}(\Phi)_{\mathbf{s}}$, since $\mathrm{Cg}_{\mathbf{A}}(\Phi)$ is a $\Sigma$-congruence and, by the 
inductive hypothesis, for every $j\in\bb{\mathbf{s}}$, $\mathrm{C}^{n}_{\mathbf{A}}(\Phi)_{s_{j}}\subseteq \mathrm{Cg}_{\mathbf{A}}(\Phi)_{s_{j}}$, we have that the pair $(\sigma^{\mathbf{A}}((x_{j})_{j\in\bb{\mathbf{s}}}, \sigma^{\mathbf{A}}((y_{j})_{j\in\bb{\mathbf{s}}})$ belongs to $\mathrm{Cg}_{\mathbf{A}}(\Phi)_{s}$.

All in all, we can affirm that $\mathrm{C}^{\omega}_{\mathbf{A}}(\Phi)\subseteq \mathrm{Cg}_{\mathbf{A}}(\Phi)$.
Hence, $\mathrm{Cg}_{\mathbf{A}}(\Phi)=\mathrm{C}^{\omega}_{\mathbf{A}}(\Phi)$.

This finishes the proof.
\end{proof}

\begin{proposition}
Let $f$ be a surjective homomorphism from the $\Sigma$-algebra $\mathbf{A}$ to the $\Sigma$-algebra 
$\mathbf{B}$ and $\Phi\subseteq A\times A$. Then 
$$
\mathrm{Cg}_{\mathbf{B}}(f^{2}[\Phi]) = f^{2}[\mathrm{Ker}(f)\vee\mathrm{Cg}_{\mathbf{A}}(\Phi)].
$$
\end{proposition}

\section{Elementary translations and translations}\label{STrans}

We next define for a $\Sigma$-algebra the concepts of elementary translation and of translation with respect to it, and provide, by using the just mentioned concepts, two characterizations of the congruences on a $\Sigma$-algebra. To this we add that the concept of translation will play an essential role in the work at hand. 

\begin{definition}\label{DETrans}
\index{translation!zerothth-order!elementary}
Let $\mathbf{A}$ be a $\Sigma$-algebra and $t\in S$. Then we will denote by $\mathrm{Etl}_{t}(\mathbf{A})$ the subset $(\mathrm{Etl}_{t}(\mathbf{A})_{s})_{s\in S}$ of $(\mathrm{Hom}(A_{t},A_{s}))_{s\in S}$ defined, for every $s\in S$, as follows: for every mapping $T\in \mathrm{Hom}(A_{t},A_{s})$, $T\in \mathrm{Etl}_{t}(\mathbf{A})_{s}$ if and only if the following condition holds
\begin{enumerate}
\item There is a word $\mathbf{s}\in  S^{\star}-\{\lambda\}$, an index $k\in \lvert \mathbf{s} \rvert$, an operation symbol $\sigma\in \Sigma_{\mathbf{s},s}$, a family $(a_{j})_{j\in i}\in\prod_{j\in k}A_{s_{j}}$, and a family $(a_{l})_{l\in \lvert \mathbf{s} \rvert-(k+1)} \in\prod_{l\in \lvert \mathbf{s} \rvert-(k+1)}A_{s_{l}}$ (recall that $k+1 = \{0, 1,\ldots,k\}$ and that $\lvert \mathbf{s} \rvert-(k+1) = \{k+1,\ldots,\lvert \mathbf{s} \rvert-1\}$) such that $s_{k} = t$ and, for every $x\in A_{t}$, 
$$T(x) =
\sigma^{\mathbf{A}}(a_{0},\ldots,a_{k-1},x,a_{k+1},\ldots,a_{\lvert \mathbf{s} \rvert-1}).$$ 
\end{enumerate}
We will sometimes use the following presentation of the elementary translations, which consists in adding an underlined space to denote where the variable will be placed
$$T =
\sigma^{\mathbf{A}}(a_{0},\ldots,a_{k-1},\underline{\quad},a_{k+1},\ldots,a_{\lvert \mathbf{s} \rvert-1}).$$ 
In this case we will say that $T$ is an \emph{elementary translation of type} $\sigma$.
We will call the elements of $\mathrm{Etl}_{t}(\mathbf{A})_{s}$ the $t$-\emph{elementary translations of sort} $s$ \emph{for} $\mathbf{A}$.
\end{definition}

\begin{definition}\label{DTrans}
\index{translation!zerothth-order}
Let $\mathbf{A}$ be a $\Sigma$-algebra and $t\in S$. Then we will denote by $\mathrm{Tl}_{t}(\mathbf{A})$ the subset $(\mathrm{Tl}_{t}(\mathbf{A})_{s})_{s\in S}$ of $(\mathrm{Hom}(A_{t},A_{s}))_{s\in S}$ defined, for every $s\in S$, as follows: for every mapping $T\in \mathrm{Hom}(A_{t},A_{s})$, $T\in \mathrm{Tl}_{t}(\mathbf{A})_{s}$ if and only if there is an $n\in \mathbb{N}-1$, a word $\mathbf{w}\in S^{n+1}$, and a family $(T_{j})_{j\in n}$ such that $w_{0} = t$, $w_{n} = s$, for every $j\in n$, $T_{j}\in \mathrm{Etl}_{w_{j}}(\mathbf{A})_{w_{j+1}}$, and $T = T_{n-1}\circ\cdots\circ T_{0}$. For translations, as for words on an alphabet, we have the notion of subtranslation of a translation, which is the counterpart of that of subword of a word. In particular, for a translation as above we will let $T'$ stand for the composition $T_{n-2}\circ \cdots \circ T_{0}$ and we will call it the \emph{maximal prefix of $T$}, and we will represent $T$ as $T_{n-1}\circ T'$ or under the form: 
$$
T=\sigma^{\mathbf{A}}\left(
a_{0},\cdots, a_{k-1}, T',
a_{k+1}, \cdots, a_{\mathbf{s}-1}
\right),
$$
where $T_{n-1}=\sigma^{\mathbf{A}}(a_{0},\ldots,a_{k-1},\underline{\quad},a_{k+1},\ldots,a_{\lvert \mathbf{s} \rvert-1})$. 
The underlined space notation can be extended to translations as well. We will say that $T$ is a \emph{translation of type $\sigma$} if the elementary translation $T_{n-1}$ is of type $\sigma$. We will call $n$ the \emph{height} of $T$ and we will denote this fact by $\bb{T}=n$. In this regard, elementary translations have height $1$ and, if $T$ is a translation of height $n$, i.e., $\bb{T}=n$, then its maximal prefix has height $n-1$, i.e., $\bb{T'}=n-1$. We will call the elements of $\mathrm{Tl}_{t}(\mathbf{A})_{s}$ the $t$-\emph{translations of sort} $s$ \emph{for} $\mathbf{A}$. In addition, for every $t\in S$, the mapping  $\mathrm{id}^{A_{t}}$ will be viewed as an element of $\mathrm{Tl}_{t}(\mathbf{A})_{t}$. The identity translation has no associated type and we will consider that it has height $0$, i.e., $\bb{\mathrm{id}^{A_{t}}}=0$. Moreover, since we consider that, for every $t\in S$, the identity mapping $\mathrm{id}^{A_{t}}$ is a traslation, we agree that an elementary translation $T\colon A_{t}\mor A_{s}$, of sort $s$ for $\mathbf{A}$, has as maximal prefix the identity at the domain.
\end{definition}

\begin{remark}
The $S\times S$-sorted set $(\mathrm{Tl}_{t}(\mathbf{A})_{s})_{(t,s)\in S\times S}$ determines a category $\mathsf{Tl}(\mathbf{A})$ whose object set is $S$ and in which, for every $(t,s)\in S\times S$, $\mathrm{Hom}_{\mathsf{Tl}(\mathbf{A})}(t,s)$, the hom-set from $t$ to $s$, is $\mathrm{Tl}_{t}(\mathbf{A})_{s}$. Therefore, for every $t\in S$, $\mathrm{End}_{\mathsf{Tl}(\mathbf{A})}(t)$ is equipped with a structure of monoid.
\end{remark}

The following are the above-mentioned characterizations of the congruences on a $\Sigma$-algebra. We note that in~\cite{m76}, on p.~199, it was announced without proof a proposition similar to that set out below.

\begin{proposition}\label{PTransCong}
Let $\mathbf{A}$ be a $\Sigma$-algebra and $\Phi$ an $S$-sorted equivalence on $A$. Then the following conditions are equivalent:
\begin{enumerate}
\item $\Phi$ is a congruence on $\mathbf{A}$.
\item $\Phi$ is closed under the elementary translations on
$\mathbf{A}$, i.e., for every every $t$, $s\in S$, every $x$, $y\in A_{t}$, and every $T\in \mathrm{Etl}_{t}(\mathbf{A})_{s}$, if $(x,y)\in \Phi_{t}$, then $(T(x),T(y))\in \Phi_{s}$.
\item $\Phi$ is closed under the translations on $\mathbf{A}$, i.e., for every every $t$, $s\in S$, every $x$, $y\in A_{t}$, and every $T\in \mathrm{Tl}_{t}(\mathbf{A})_{s}$, if $(x,y)\in \Phi_{t}$, then $(T(x),T(y))\in \Phi_{s}$.
\end{enumerate}
\end{proposition}

\begin{proof}
Let us first prove that (1) and (2) are equivalent.

Let us suppose that $\Phi$ is a congruence on $\mathbf{A}$. We want to show that $\Phi$ is closed under the elementary translations on $\mathbf{A}$. Let $t$ and $s$ be elements of $S$ and $T$ a $t$-elementary translation of sort $s$ for $\mathbf{A}$. Then $T\colon A_{t}\mor A_{s}$ and there is a word $\mathbf{s}\in S^{\star}-\{\lambda\}$, an index $k\in \lvert \mathbf{s} \rvert$, an operation symbol $\sigma\in
\Sigma_{\mathbf{s},s}$, a family $(a_{j})_{j\in k}\in\prod_{j\in k}A_{s_{j}}$, and a family $(a_{l})_{l\in \lvert \mathbf{s} \rvert-(k+1)}
\in\prod_{kl\in \lvert \mathbf{s} \rvert-(k+1)}A_{s_{l}}$ such that $s_{k} = t$ and, for every $z\in A_{t}$, $T(z) =
\sigma^{\mathbf{A}}(a_{0},\ldots,a_{k-1},z,a_{k+1},\ldots,a_{\lvert \mathbf{s} \rvert-1})$. Let $x$ and $y$ be elements of $A_{t}$ such that $(x,y)\in \Phi_{t}$. Since, for every $j\in k$, $(a_{j},a_{j})\in \Phi_{s_{j}}$, for every $l\in \lvert \mathbf{s} \rvert-(k+1)$, $(a_{l},a_{l})\in \Phi_{s_{l}}$, and, in addition, $(x,y)\in \Phi_{t} = \Phi_{s_{k}}$, then $(T(x),T(y))\in \Phi_{s}$.

Reciprocally, let us suppose that, for every $t$, $s\in S$, every $x$, $y\in A_{t}$, and every $T\in \mathrm{Etl}_{t}(\mathbf{A})_{s}$, if $(x,y)\in \Phi_{t}$, then $(T(x),T(y))\in \Phi_{s}$. We want to show that  $\Phi$ is a congruence on $\mathbf{A}$. Let $(\mathbf{s},s)\in (S^{\star}-\{\lambda\})\times S$, $\sigma\in\Sigma_{\mathbf{s},s}$,
and $(a_{i})_{i\in\lvert \mathbf{s} \rvert}$ and $ (b_{i})_{i\in\lvert \mathbf{s} \rvert}$ be families in $ A_{\mathbf{s}}$ such that, for every $i\in \lvert \mathbf{s} \rvert$ we have that $(a_{i}, b_{i})\in \Phi_{s_{i}}$. We now define, for every $i\in\lvert \mathbf{s} \rvert$, $T_{i}$, the $s_{i}$-elementary translation of sort $s$ for $\mathbf{A}$, as the mapping from $A_{s_{i}}$ to $A_{s}$ which sends $x\in A_{s_{i}}$ to $\sigma^{\mathbf{A}}(b_{0},\ldots,b_{i-1},x,a_{i+1},\ldots,a_{\lvert \mathbf{s} \rvert-1}) \in A_{s}$. Then $\sigma^{\mathbf{A}}(a_{0},\ldots,a_{\lvert \mathbf{s} \rvert-1}) = T_{0}(a_{0})$ and $(T_{0}(a_{0}),T_{0}(b_{0}))\in \Phi_{s_{0}}$. But $T_{0}(b_{0}) = T_{1}(a_{1})$ and $(T_{1}(a_{1}),T_{1}(b_{1}))\in \Phi_{s_{1}}$. By proceeding in the same way we, finally, come to $T_{\lvert \mathbf{s} \rvert-2}(b_{\lvert \mathbf{s} \rvert-2}) = T_{\lvert \mathbf{s} \rvert-1}(a_{\lvert \mathbf{s} \rvert-1})$, $(T_{\lvert \mathbf{s} \rvert-1}(a_{\lvert \mathbf{s} \rvert-1}), T_{\lvert \mathbf{s} \rvert-1}(b_{\lvert \mathbf{s} \rvert-1}))\in \Phi_{s_{\lvert \mathbf{s} \rvert-1}}$, and $T_{\lvert \mathbf{s} \rvert-1}(b_{\lvert \mathbf{s} \rvert-1})\! =\! \sigma^{\mathbf{A}}(b_{0},\ldots,b_{\lvert \mathbf{s} \rvert-1})$. Therefore $(\sigma^{\mathbf{A}}((a_{j})_{j\in\bb{\mathbf{s}}}), \sigma^{\mathbf{A}}((b_{j})_{j\in\bb{\mathbf{s}}}))\in \Phi_{s}$.

We shall now proceed to verify that (2) and (3) are equivalent.

Since every elementary translations on $\mathbf{A}$ is a translation on $\mathbf{A}$, it is obvious that if $\Phi$ is closed under the translations on $\mathbf{A}$, then $\Phi$ is closed under the elementary translations on $\mathbf{A}$.

Reciprocally, let us suppose that $\Phi$ is closed under the elementary translations on $\mathbf{A}$. We want to show that $\Phi$ is closed under the translations on $\mathbf{A}$. Let $t$ and $s$ be elements of $S$, $x$, $y$ elements of $A_{t}$, $T\in \mathrm{Tl}_{t}(\mathbf{A})_{s}$, and let us suppose that $(x,y)\in \Phi_{t}$. Then there is an $n\in \mathbb{N}-1$, a word $\mathbf{w}\in S^{n+1}$, and a family $(T_{j})_{j\in n}$ such that $w_{0} = t$, $w_{n} = s$, for every $j\in n$, 
$T_{j}\in \mathrm{Etl}_{w_{j}}(\mathbf{A})_{w_{j+1}}$,
and $T = T_{n-1}\circ\cdots\circ T_{0}$. Then, from $(x,y)\in \Phi_{t} = \Phi_{w_{0}}$, we infer that $(T_{0}(x),T_{0}(y))\in \Phi_{w_{1}}$. By proceeding in the same way we, finally, come to $(T_{n-1}(\ldots(T_{0}(x))\ldots),T_{n-1}(\ldots(T_{0}(y))\ldots))\in \Phi_{s} = \Phi_{w_{n}}$, i.e., to $(T(x),T(y))\in \Phi_{s}$.
\end{proof}


We next investigate the relationships between the translations and the homomorphisms between $\Sigma$-algebras. The first result is a direct consequence of Proposition~\ref{PTransCong}.

\begin{corollary}\label{CTransCong} Let $\mathbf{A}$ and $\mathbf{B}$ be two $\Sigma$-algebras and $f\colon \mathbf{A}\mor \mathbf{B}$ a $\Sigma$-homomorphism. For $s,t\in S$, let $T$ be a translation on $A$ in $\mathrm{Tl}_{t}(\mathbf{A})_{s}$. Let $x,y\in A_{t}$ be such that $f_{t}(x)=f_{t}(y)$, then $f_{s}(T(x))=f_{s}(T(y))$.
\end{corollary}
\begin{proof}
It follows directly from Proposition~\ref{PTransCong} and from the fact that 	$\mathrm{Ker}(f)$ is a $\Sigma$-congruence on $\mathbf{A}$.
\end{proof}

Now we elucidate the relationships that exist between translations and homomorphisms.

\begin{proposition}\label{PTransHom}
Let $\mathbf{A}$ and $\mathbf{B}$ be two $\Sigma$-algebras and $f\colon \mathbf{A}\mor \mathbf{B}$ a $\Sigma$-homomorphism. Then, for every $t$, $s\in S$ and every $T\in \mathrm{Tl}_{t}(\mathbf{A})_{s}$, there exists a $T^{f}\in\mathrm{Tl}_{t}(\mathbf{B})_{s}$ such that $f_{s}\circ T = T^{f}\circ f_{t}$. Moreover, if $f$ is an epimorphism, then, for every $t$, $s\in S$ and every $U\in\mathrm{Tl}_{t}(\mathbf{B})_{s}$, there exists a $T\in \mathrm{Tl}_{t}(\mathbf{A})_{s}$ such that $T^{f} = U$.
\end{proposition}

\begin{proof}
If $T$ is a $t$-elementary translation of sort $s$ for $\mathbf{A}$, for some $t$ and $s$ in $S$, then there is a word $\mathbf{s}\in S^{\star}-\{\lambda\}$, an index $k\in \lvert \mathbf{s} \rvert$, an operation symbol $\sigma\in \Sigma_{\mathbf{s},s}$, a family $(a_{j})_{j\in k}\in\prod_{j\in k}A_{s_{j}}$, and a family $(a_{l})_{l\in \lvert \mathbf{s} \rvert-(k+1)}\in\prod_{l\in \lvert \mathbf{s} \rvert-(k+1)}A_{s_{l}}$ such that $s_{k} = t$ and, for every $x\in A_{t}$, $T(x) = \sigma^{\mathbf{A}}(a_{0},\ldots,a_{k-1},x,a_{k+1},\ldots,a_{\lvert \mathbf{s} \rvert-1})$. Then it suffices to take as $T^{f}$ precisely the mapping from $B_{t}$ to $B_{s}$ defined, for every $y\in B_{t}$, as follows:
$$
T^{f}(y) = \sigma^{\mathbf{B}}(f_{s_{0}}(a_{0}),\ldots,f_{s_{k-1}}(a_{k-1}),y,f_{s_{k+1}}(a_{k+1}),\ldots,f_{s_{\lvert \mathbf{s} \rvert-1}}(a_{\lvert \mathbf{s} \rvert-1})).
$$

If $T\in\mathrm{Tl}_{t}(\mathbf{A})_{s}$, for some $t$ and $s$ in $S$, and $T$ is not a $t$-elementary translation of sort $s$ for $\mathbf{A}$, then there is an $n\in \mathbb{N}-1$, a word $\mathbf{w}\in S^{n+1}$, and a family $(T_{j})_{j\in n}$ such that $w_{0} = t$, $w_{n} = s$, for every $j\in n$, 
$T_{j}\in \mathrm{Etl}_{w_{j}}(\mathbf{A})_{w_{j+1}}$
and $T = T_{n-1}\circ\cdots\circ T_{0}$. Then it suffices to take as $T^{f}$ precisely the mapping from $B_{t}$ to $B_{s}$ defined as
$T^{f} = T_{n-1}^{f}\circ\cdots\circ T_{0}^{f}$. Let us note that if $T$ is $\mathrm{id}^{{A}_{t}}$, for some $t\in S$, then it suffices to take as $T^{f}$ precisely $\mathrm{id}^{{B}_{t}}$.

We next prove that if $f$ is an epimorphism, then, for every $t$, $s\in S$ and every $U\in\mathrm{Tl}_{t}(\mathbf{B})_{s}$, there exists a $T\in \mathrm{Tl}_{t}(\mathbf{A})_{s}$ such that $T^{f} = U$. If $U$ is a $t$-elementary translation of sort $s$ for $\mathbf{B}$, for some $t$ and $s$ in $S$, then and there is a word $\mathbf{s}\in S^{\star}-\{\lambda\}$, an index $k\in \lvert \mathbf{s} \rvert$, an operation symbol $\sigma\in\Sigma_{\mathbf{s},s}$, a family $(b_{j})_{j\in k}\in\prod_{j\in k}B_{s_{j}}$, and a family $(b_{l})_{l\in \lvert \mathbf{s} \rvert-(k+1)}\in\prod_{l\in \lvert \mathbf{s} \rvert-(k+1)}B_{s_{l}}$ such that $s_{k} = t$ and, for every $y\in B_{t}$, $U(y) = \sigma^{\mathbf{B}}(b_{0},\ldots,b_{k-1},y,b_{k+1},\ldots,b_{\lvert \mathbf{s} \rvert-1})$. Then, since $f$ is an epimorphism, there exists a family $(a_{j})_{j\in k}\in\prod_{j\in k}A_{s_{j}}$ and a family $(a_{l})_{l\in \lvert \mathbf{s} \rvert-(l+1)}\in\prod_{l\in \lvert \mathbf{s} \rvert-(k+1)}A_{s_{l}}$ such that, for every $j\in k$, $f_{s_{j}}(a_{j}) = b_{j}$, and, for every $l\in \lvert \mathbf{s} \rvert-(k+1)$, $f_{s_{l}}(a_{l}) = b_{l}$. Then, after fixing $(a_{j})_{j\in k}$ and $(a_{l})_{l\in \lvert \mathbf{s} \rvert-(k+1)}$, it suffices to take as $T$ precisely the mapping from $A_{t}$ to $A_{s}$ defined, for every $x\in A_{t}$, as follows:
$$
T(x) = \sigma^{\mathbf{A}}(a_{0},\ldots,a_{k-1},x,a_{k+1},\ldots,a_{\lvert \mathbf{s} \rvert-1}).
$$

If $U\in\mathrm{Tl}_{t}(\mathbf{B})_{s}$, for some $t$ and $s$ in $S$, and $U$ is not a $t$-elementary translation of sort $s$ for $\mathbf{B}$, then there is an $n\in \mathbb{N}-1$, a word $\mathbf{w}\in S^{n+1}$, and a family $(U_{j})_{j\in n}$ such that $w_{0} = t$, $w_{n} = s$, for every $j\in n$
$U_{j}\in \mathrm{Etl}_{w_{j}}(\mathbf{A})_{w_{j+1}}$, 
and $U = U_{n-1}\circ\cdots\circ U_{0}$. Then, after choosing, for every $i\in n$, a $T_{i}$ such that $T_{i}^{f} = U_{i}$, it suffices to take as $T$ precisely the mapping from $A_{t}$ to $A_{s}$ defined as $T = T_{n-1}\circ\cdots\circ T_{0}$.
\end{proof}

\begin{corollary}\label{PTransHomComp}
Let $\mathbf{A}$ and $\mathbf{B}$ be two $\Sigma$-algebras, $f\colon \mathbf{A}\mor \mathbf{B}$ a $\Sigma$-ho\-mo\-mor\-phism, $t$, $s\in S$, $T\in \mathrm{Tl}_{t}(\mathbf{A})_{s}$, $\mathbf{w}\in S^{\star}$ and $\sigma\in \Sigma_{\mathbf{w},t}$. Then $f_{s}\circ T\circ \sigma^{\mathbf{A}} = T^{f}\circ \sigma^{\mathbf{B}}\circ f_{\mathbf{w}}$.
\end{corollary}

\section{
\texorpdfstring
{The free $\Sigma$-algebra on an $S$-sorted set}
{The free algebra on a many-sorted set}
}\label{SAlgFree}

We next state that the forgetful functor $\mathrm{G}_{\Sigma}$ from $\mathsf{Alg}(\Sigma)$ to
$\mathsf{Set}^{S}$ has a left adjoint $\mathbf{T}_{\Sigma}$ which assigns to an $S$-sorted set $X$ the free $\Sigma$-algebra $\mathbf{T}_{\Sigma}(X)$ on $X$. We later state the universal property of the free many-sorted algebra and provide a characterization by means of the notion of Dedekind Peano-algebras. Furthermore, for free many-sorted algebras, we establish that the subelement preorder on them is an Artinian order and, for a term in a free many-sorted algebra, we define the sorted set of the subterms of it.

Let us note that in what follows, to construct the algebra of $\Sigma$-rows in $X$, and the free $\Sigma$-algebra on $X$, since neither the $S$-sorted signature $\Sigma$ nor the $S$-sorted set $X$ are subject to any constraint, coproducts must necessarily be used.

\begin{definition}
Let $X$ be an $S$-sorted set. The \emph{algebra of} $\Sigma$-\emph{rows in} $X$, denoted by $\mathbf{W}_{\Sigma}(X)$, is defined as follows:
\begin{enumerate}
\item The underlying $S$-sorted set of $\mathbf{W}_{\Sigma}(X)$, written as $\mathrm{W}_{\Sigma}(X)$, is precisely the $S$-sorted set $((\coprod\Sigma \amalg \coprod X)^{\star})_{s\in S}$, i.e., the mapping from $S$ to $\boldsymbol{\mathcal{U}}$ constantly $(\coprod\Sigma \amalg \coprod X)^{\star}$, where $(\coprod\Sigma \amalg \coprod X)^{\star}$ is the set of all words on the set $\coprod\Sigma \amalg \coprod X$, i.e., on the set
      $$
      \textstyle
      [(\bigcup_{(\mathbf{s},s)\in S^{\star}\times S}(\Sigma_{\mathbf{s},s}\times\{(\mathbf{s},s)\}))\times\{0\}]\cup
      [(\bigcup_{s\in S}(X_{s}\times\{s\}))\times\{1\}].
      $$

\item For every $(\mathbf{s},s)\in S^{\star}\times S$, and every $\sigma\in\Sigma_{\mathbf{s},s}$, the structural operation $F_{\sigma}$ associated to $\sigma$ is the mapping from $\mathrm{W}_{\Sigma}(X)_{\mathbf{s}}$ to ${\mathrm{W}_{\Sigma}(X)}_{s}$ which sends $(P_{j})_{j\in\bb{\mathbf{s}}} \in \mathrm{W}_{\Sigma}(X)_{\mathbf{s}}$ to $(\sigma)\curlywedge\concat_{j\in\bb{\mathbf{s}}}P_{j} \in {\mathrm{W}_{\Sigma}(X)}_{s}$, where, for every $(\mathbf{s},s)\in S^{\star}\times S$, and every $\sigma\in\Sigma_{\mathbf{s},s}$,
      $(\sigma)$ stands for $(((\sigma,(\mathbf{s},s)),0))$, which is the value at $\sigma$ of the canonical mapping from $\Sigma_{\mathbf{s},s}$ to $(\coprod\Sigma \amalg \coprod X)^{\star}$.
\end{enumerate}
\end{definition}

\begin{definition}\label{DAlgFree}
The \emph{free} $\Sigma$-\emph{algebra on} an $S$-sorted set $X$, denoted by $\mathbf{T}_{\Sigma}(X)$, is the $\Sigma$-algebra determined by $\mathrm{Sg}_{\mathbf{W}_{\Sigma}(X)}((\{(x)\mid x\in X_{s}\})_{s\in S})$, the subalgebra of $\mathbf{W}_{\Sigma}(X)$ generated by $(\{(x)\mid x\in X_{s}\})_{s\in S}$, where, for every $s\in S$ and every $x\in X_{s}$, $(x)$ stands for $(((x,s),1))$, which is the value at $x$ of the canonical mapping from $X_{s}$ to $(\coprod\Sigma \amalg \coprod X)^{\star}$.
We will denote by $\mathrm{T}_{\Sigma}(X)$ the underlying $S$\nobreakdash-sorted set of $\mathbf{T}_{\Sigma}(X)$ and, for $s\in S$, we will call the elements of $\mathrm{T}_{\Sigma}(X)_{s}$ \emph{terms of type} $s$ \emph{with variables in} $X$ or  $(X,s)$-\emph{terms}.
\end{definition}

\begin{remark}
Since $(\{(x)\mid x\in X_{s}\})_{s\in S}$ is a generating subset of $\mathbf{T}_{\Sigma}(X)$, to prove that a subset $\mathcal{T}$ of $\mathrm{T}_{\Sigma}(X)$ is equal to $\mathrm{T}_{\Sigma}(X)$ it suffices, by Proposition~\ref{PPAI}, to show: (1) $(\{(x)\mid x\in X_{s}\})_{s\in S}\subseteq \mathcal{T}$ (algebraic induction basis) and (2) $\mathcal{T}$ is a subalgebra of $\mathbf{T}_{\Sigma}(X)$ (algebraic induction step).
\end{remark}

In the many-sorted case we have, as in the single-sorted case, the following characterization of the elements of $\mathrm{T}_{\Sigma}(X)_{s}$, for $s\in S$.

\begin{proposition}\label{PTermChar}
Let $X$ be an $S$-sorted set. Then, for every $s\in S$ and every $P\in
\mathrm{T}_{\Sigma}(X)_{s}$, we have that $P$ is a term of type $s$ with variables in $X$ if and only if $P = (x)$, for a unique $x\in X_{s}$, or $P = (\sigma)$, for a unique $\sigma\in\Sigma_{\lambda,s}$, or $P = (\sigma)\curlywedge\concat(P_{j})_{j\bb{\mathbf{s}}}$, for a unique $\mathbf{s}\in S^{\star}-\{\lambda\}$, a unique $\sigma\in\Sigma_{\mathbf{s},s}$, and a unique family $(P_{j})_{j\in\bb{\mathbf{s}}}\in\mathrm{T}_{\Sigma}(X)_{\mathbf{s}}$.

Moreover, the three possibilities are mutually exclusive. From now on, for simplicity of notation, we will write $x$, $\sigma^{\mathbf{T}_{\Sigma}(X)}$, and $\sigma^{\mathbf{T}_{\Sigma}(X)}(P_{0},\ldots,P_{\bb{\mathbf{s}}-1})$ or $\sigma^{\mathbf{T}_{\Sigma}(X)}((P_{j})_{j\in \bb{\mathbf{s}}})$ instead of $(x)$, $(\sigma)$, and $(\sigma)\curlywedge\concat(P_{j})_{j\in\bb{\mathbf{s}}}$, respectively. Thus, in particular, the structural operation $F_{\sigma}$ (or more accurately $F_{\sigma}^{\mathbf{T}_{\Sigma}(X)}$) associated with $\sigma$ is identified with $\sigma^{\mathbf{T}_{\Sigma}(X)}$.
\end{proposition}

Since the terms of the last class of the proposition just stated, as well as the first letter of them, will be constantly used in this work, in the following definition, we give them special names.

\begin{definition}
If a term $P\in \mathrm{T}_{\Sigma}(X)_{s}$ of type $(X,s)$ is of the form $$\sigma^{\mathbf{T}_{\Sigma}(X)}((P_{j})_{j\in \bb{\mathbf{s}}}),$$ for a unique $\mathbf{s}\in S^{\star}-\{\lambda\}$, a unique $\sigma\in\Sigma_{\mathbf{s},s}$, and a unique family $(P_{j})_{j\in\bb{\mathbf{s}}}\in\mathrm{T}_{\Sigma}(X)_{\mathbf{s}}$, then it will be called a \emph{complex term} and the first letter of it (which is the operation symbol $\sigma$) will be called the \emph{head} of $P$ and will be denoted by $\mathrm{hd}(P)$.
\end{definition}

We next distinguish, from among the non-empty finite families of terms of the same type with variables in a given $S$-sorted set, those that are  head-constant. And we do such a thing because in this work some specific types of  head-constant families of terms will play an important role. 

\begin{definition} \label{DHeadCt}
Let $X$ be an $S$-sorted set, $s\in S$, $n\in \mathbb{N}-1$, and $(P_{i})_{i\in n}\in \mathrm{T}_{\Sigma}(X)_{s}^{n}$. If, for every $i\in n$, $P_{i}$ is a complex term and, for every $i$, $j\in n$, $\mathrm{hd}(P_{i}) = \mathrm{hd}(P_{j})$, then we will say that the family $(P_{i})_{i\in n}$ is \emph{head-constant}.  
Therefore, $(P_{i})_{i\in n}$ is complex and head-constant if, and only if, there exists a unique word $\mathbf{s}\in S^{\star}-\{\lambda\}$, a unique operation
symbol $\sigma\in \Sigma_{\mathbf{s},s}$, and a unique family of path terms $((P_{i,j})_{j\in\bb{\mathbf{s}}})_{i\in n}\in \mathrm{T}_{\Sigma}(X)_{\mathbf{s}}^{n}$ such that, for every $i\in n$,
$$
P_{i}=\sigma^{\mathbf{T}_{\Sigma}(X)}((P_{i,j})_{j\in\bb{\mathbf{s}}}).
$$ 
\end{definition}


From Proposition~\ref{PTermChar} follows, immediately, the universal property of the free $\Sigma$-algebra on an $S$-sorted set $X$, as stated in the subsequent proposition.

\begin{proposition}\label{PPropUniv}
For every $S$-sorted set $X$, the pair $(\eta^{X},\mathbf{T}_{\Sigma}(X))$, where $\eta^{X}$, the
\emph{insertion of (the $S$-sorted set of generators)} $X$ \emph{into} $\mathrm{T}_{\Sigma}(X)$, is the co-restric\-tion to
$\mathrm{T}_{\Sigma}(X)$ of the canonical embedding of $X$ into $\mathrm{W}_{\Sigma}(X)$, has the following universal property: for every $\Sigma$-algebra $\mathbf{A}$ and every $S$-sorted mapping $f\colon X\mor A$, there exists a unique homomorphism $f^{\sharp}\colon\mathbf{T}_{\Sigma}(X)\mor\mathbf{A}$ such that $f^{\sharp}\circ \eta^{X} = f$.
\end{proposition}

\begin{proof}
For every $s\in S$ and every $(X,s)$-term $P$, the $s$-th coordinate $f^{\sharp}_{s}$ of $f^{\sharp}$ is defined recursively as follows: $f^{\sharp}_{s}(x) = f_{s}(x)$, if $P = x$; $f^{\sharp}_{s}(\sigma^{\mathbf{T}_{\Sigma}(X)}) = \sigma^{\mathbf{A}}$, if $P = \sigma^{\mathbf{T}_{\Sigma}(X)}$; and, finally,  $f^{\sharp}_{s}(\sigma^{\mathbf{T}_{\Sigma}(X)}((P_{j})_{j\in\bb{\mathbf{s}}})) = \sigma^{\mathbf{A}}((f^{\sharp}_{s_{j}}(P_{j}))_{j\in\bb{\mathbf{s}}})$, if $P = \sigma^{\mathbf{T}_{\Sigma}(X)}((P_{j})_{j\in\bb{\mathbf{s}}})$.
\end{proof}

The just stated proposition allows us to carry out definitions by algebraic recursion on a free many-sorted algebra as indeed we will be doing throughout  this paper.


\begin{corollary} \label{CTermAdj}
The functor $\mathbf{T}_{\Sigma}$, which sends an $S$-sorted set $X$ to $\mathbf{T}_{\Sigma}(X)$ and an $S$-sorted mapping $f$ from $X$ to $Y$ to $f^{@} (= (\eta^{Y}\circ f)^{\sharp})$, the unique homomorphism from $\mathbf{T}_{\Sigma}(X)$ to $\mathbf{T}_{\Sigma}(Y)$ such that $f^{@}\circ\eta^{X} = \eta^{Y}\circ f$, is left adjoint for the forgetful functor $\mathrm{G}_{\Sigma}$ from $\mathsf{Alg}(\Sigma)$ to $\mathsf{Set}^{S}$.
\end{corollary}

\begin{remark}
Every many-sorted \emph{algebra} in a variety is isomorphic to a quotient of a free algebra in the variety. 
\end{remark}

It is possible to give an internal characterization of the free algebras by means of the Dedekind-Peano algebras.

\begin{definition}
Let $\mathbf{A}$ be a $\Sigma$-algebra.  We will say that $\mathbf{A}$ is a \emph{Dedekind-Peano $\Sigma$-algebra}, abbreviated to DP $\Sigma$-\emph{algebra} when this is unlikely to cause confusion, if the following axioms hold
\begin{enumerate}
\item[DP1.] For every $(\mathbf{s},s)\in S^{\star}\times S$ and every $\sigma\in\Sigma_{\mathbf{s},s}$, $\sigma^{\mathbf{A}}\colon A_{\mathbf{s}}\mor A_{s}$ is injective.

\item[DP2.] For every $s\in S $ and every $\sigma, \tau \in\Sigma_{\cdot,s}$, if $\sigma\neq\tau$, then $\mathrm{Im}(\sigma^{\mathbf{A}})\cap \mathrm{Im}(\tau^{\mathbf{A}}) = \varnothing$.

\item[DP3.] $\mathrm{Sg}_{\mathbf{A}}(A-(\bigcup_{\sigma\in\Sigma_{\cdot,s}}\mathrm{Im}(\sigma^{\mathbf{A}})_{s\in S}) = A$.
\end{enumerate}
We call the $S$-sorted set $A-(\bigcup_{\sigma\in\Sigma_{\cdot,s}}\mathrm{Im}(\sigma^{\mathbf{A}})_{s\in S})$
the \emph{basis of Dedekind-Peano of} $\mathbf{A}$, and we denote it by $\mathrm{B}(\mathbf{A})$.
\end{definition}

We remark that the axioms DP1 and DP2 are equivalent to the following axiom
\begin{enumerate}
\item[DP4.] For every $s\in S$, $\sigma, \tau \in\Sigma_{\cdot,s}$, $a\in A_{\mathrm{ar}(\sigma)}$, $b\in
A_{\mathrm{ar}(\tau)}$, if $\sigma^{\mathbf{A}}(a)=\tau^{\mathbf{A}}(b)$ then $\sigma=\tau$ and $a=b$.
\end{enumerate}

Moreover, if $\mathbf{A}$ is a DP $\Sigma$-algebra, then $\mathrm{B}(\mathbf{A})$ is $\bigcap\{X\subseteq A\mid \mathrm{Sg}_{\mathbf{A}}(X)=A\}$

\begin{proposition}
Let $\mathbf{A}$ be a DP $\Sigma$-algebra. Then $\mathbf{A}$ is isomorphic
to $\mathbf{T}_{\Sigma}(\mathrm{B}(\mathbf{A}))$.
\end{proposition}

\begin{proposition}
Let $X$ be an $S$-sorted set. Then $\mathbf{T}_{\Sigma}(X)$ is a DP $\Sigma$-algebra.
\end{proposition}

\begin{proposition}
Let $\mathbf{A}$ be a DP $\Sigma$-algebra. Then $\leq_{\mathbf{A}}$ is antisymmetric and does not have strictly descending $\omega_{0}$-chains, i.e., is an Artinian order.
\end{proposition}

\begin{corollary}
Let $X$ be an $S$-sorted set. Then $\leq_{\mathbf{T}_{\Sigma}(X)}$ is antisymmetric and does not have strictly descending $\omega_{0}$-chains, i.e., is an Artinian order.
\end{corollary}

\begin{remark}\label{RTermOrd}
Let us recall from Definition~\ref{DAlgOrd} that $\prec_{\mathbf{T}_{\Sigma}(X)}$ denotes the binary relation on $\coprod \mathrm{T}_{\Sigma}(X)$ consisting of the ordered pairs $((Q,t),(P,s))\in(\coprod \mathrm{T}_{\Sigma}(X))^{2}$ for which there exists a word $\mathbf{s}\in S^{\star}-\{\lambda\}$, an operation symbol $\sigma\in \Sigma_{\mathbf{s},s}$ and a family of terms $(P_{j})_{j\in\bb{\mathbf{s}}}\in \mathrm{T}_{\Sigma}(X)_{\mathbf{s}}$ such that $P=\sigma^{\mathbf{T}_{\Sigma}(X)}((P_{j})_{j\in\bb{\mathbf{s}}})$ and, for some $j\in\bb{\mathbf{s}}$, $s_{j}=t$ and $Q=P_{j}$. 

We denote by $\leq_{\mathbf{T}_{\Sigma}(X)}$ the reflexive and transitive closure of $\prec_{\mathbf{T}_{\Sigma}(S)}$, i.e., the preorder on $\coprod \mathrm{T}_{\Sigma}(X)$  generated by $\prec_{\mathbf{T}_{\Sigma}(X)}$, and by $<_{\mathbf{T}_{\Sigma}(X)}$ the transitive closure of $\prec_{\mathbf{T}_{\Sigma}(X)}$.

Note that the binary relation $\prec_{\mathbf{T}_{\Sigma}(X)}$ on $\coprod \mathrm{T}_{\Sigma}(X)$ is irreflexive.

The preorder $\leq_{\mathbf{T}_{\Sigma}(X)}$ on $\coprod \mathrm{T}_{\Sigma}(X)$ is $\bigcup_{n\in\mathbb{N}}\prec_{\mathbf{T}_{\Sigma}(X)}^{n}$, where $\prec_{\mathbf{T}_{\Sigma}(X)}^{0}$ is the diagonal of $\coprod \mathrm{T}_{\Sigma}(X)$ and, for $n\in\mathbb{N}$, 
$$\prec_{\mathbf{T}_{\Sigma}(X)}^{n+1} = \prec_{\mathbf{T}_{\Sigma}(X)}^{n}\circ \prec_{\mathbf{T}_{\Sigma}(X)}.$$ 

Thus, for every $((Q,t),(P,s))\in (\mathbf{T}_{\Sigma}(X))^{2}$, $(Q,t)\leq_{\mathbf{T}_{\Sigma}(X)}(P,s)$ if and only if $s = t$ and $Q=P$ or there exists  a natural number $m\in\mathbb{N}-\{0\}$, a word $\mathbf{c}\in S^{\star}$ of length $\bb{\mathbf{c}}=m+1$, and a family of terms $(R_{k})_{k\in\bb{\mathbf{c}}}$ in $\mathrm{T}_{\Sigma}(X)_{\mathbf{c}}$, such that $c_{0}=t$,  $R_{0}=Q$, $c_{m}=s$, $R_{m}=P$ and, for every $k\in m$, $(R_{k},c_{k})\prec_{\mathrm{T}_{\Sigma}(X)}(R_{k+1},c_{k+1})$. Moreover, 
$<_{\mathbf{T}_{\Sigma}(X)}$, the transitive closure of $\prec_{\mathbf{T}_{\Sigma}(X)}$, is 
$\bigcup_{n\in\mathbb{N}-\{0\}}\prec_{\mathbf{T}_{\Sigma}(X)}^{n}$.
\end{remark}


\begin{definition}
Let $X$ be an $S$-sorted set, $s\in S$, and $P\in\mathrm{T}_{\Sigma}(X)_{s}$. Then the $S$-sorted set of all \emph{subterms of} $P$, denoted by $\mathrm{Subt}_{\mathbf{T}_{\Sigma}(X)}(P)$ or, simply, by $\mathrm{Subt}_{\Sigma}(P)$ when no confusion can arise, is defined as follows:
$$
\mathrm{Subt}_{\Sigma}(P) = (\{Q\in \mathrm{T}_{\Sigma}(X)_{t} \mid (Q,t)
\leq_{\mathbf{T}_{\Sigma}(X)} (P,s)\})_{t \in S }.
$$
\end{definition}

\begin{remark}
$\mathrm{Subt}_{\Sigma}(P)\in \mathrm{Sub}_{\mathrm{f}}(\mathrm{T}_{\Sigma}(X))$. Moreover, $\mathrm{Subt}_{\Sigma}(P)$ can also be characterized as the smallest subset $\mathcal{L}$ of $\mathrm{T}_{\Sigma}(X)$ which satisfies the following conditions: 
\begin{enumerate}
\item $P\in \mathcal{L}_{s}$; and
\item for every $(\mathbf{s},s)\in S^{\star}\times S$, every $\sigma\in \Sigma_{\mathbf{s},s}$, and every $(Q_{j})_{j\in\bb{\mathbf{s}}}\in \mathrm{T}_{\Sigma}(X)_{\mathbf{s}}$, if $\sigma^{\mathbf{T}_{\Sigma}(X)}((Q_{j})_{j\in\bb{\mathbf{s}}})\in \mathcal{L}_{s}$, then, for every $j\in\bb{\mathbf{s}}$, $Q_{j}\in \mathcal{L}_{s_{j}}$.
\end{enumerate}
Note that the second condition is exactly the converse of the defining condition of the concept of subalgebra of a $\Sigma$-algebra.
\end{remark}

\begin{definition} Let $X$ be an $S$-sorted set and let $\mathbf{T}_{\Sigma}(X)$ be the free many-sorted $\Sigma$-algebra on $X$. For a natural number $n\in\mathbb{N}$ we will use the shorthand $\mathrm{B}^{n}_{\Sigma}(X)$ for the $n$-th Borel set of $X$ in $\mathbf{T}_{\Sigma}(X)$, that is 
$$
\mathrm{B}^{n}_{\Sigma}(X)=\mathrm{B}^{n}_{\mathbf{T}_{\Sigma}(X)}(X).
$$
Let us note that the elements of $\mathrm{T}_{\Sigma}(X)-\mathrm{B}^{0}_{\Sigma}(X)$ are precisely the complex terms.
\end{definition}

Following this we associate to every term for $\Sigma$ of type $(X,s)$ its $S$-sorted set of variables. 

\begin{definition}
Let $X$ be an $S$-sorted set.  Then $\mathbf{Fin}(X)$ is the $\Sigma$-algebra which has as underlying
$S$-sorted set $\mathrm{Fin}^{S}(X) = (\mathrm{Sub}_{\mathrm{f}}(X))_{s\in S}$, i.e., the $S$-sorted set constantly $\mathrm{Sub}_{\mathrm{f}}(X)$, and, for every $(\mathbf{s},s)\in S^{\star}\times S$ and every $\sigma\in \Sigma_{\mathbf{s},s}$, as operation $\sigma^{\mathbf{Fin}(X)}\colon \mathrm{Sub}_{\mathrm{f}}(X)^{\bb{\mathbf{s}}}\mor
\mathrm{Sub}_{\mathrm{f}}(X)$ that one defined, for every $(K^{j})_{j\in\bb{\mathbf{s}}}\in \mathrm{Sub}_{\mathrm{f}}(X)^{\bb{\mathbf{s}}}$, as $\sigma^{\mathbf{Fin}(X)}((K^{j})_{j\in\bb{\mathbf{s}}})= \bigcup_{j\in\bb{\mathbf{s}}}K^{j}$. Let
$\delta^{S,X}=(\delta^{S,X}_{s})_{s\in S}$ be the $S$-sorted mapping from $X$ to $\mathrm{Fin}^{S}(X)$ defined in Remark~\ref{RDelta}. 
Then we will denote by $\mathrm{Var}^{S,X}$ the unique homomorphism $\ext{(\delta^{S,X})}$ from $\mathbf{T}_{\Sigma}(X)$ to
$\mathbf{Fin}(X)$ such that the diagram in Figure~\ref{FVar} commutes. 
For a sort $s\in S$ and a term $P\in \mathrm{T}_{\Sigma}(X)_{s}$ we will call $\mathrm{Var}^{S,X}_{s}(P) (\in \mathrm{Sub}_{\mathrm{f}}(X))$ the $S$-sorted set of \emph{variables of} $P$ (or \emph{occurring in} $P$). Moreover, when this is unlikely to cause confusion, we will write $\mathrm{Var}^{X}(P)$ or, simply, $\mathrm{Var}(P)$ for $\mathrm{Var}^{S,X}_{s}(P)$.
\end{definition}

\begin{remark}
Let $s$ be a sort in $S$ and $P\in \mathrm{T}_{\Sigma}(X)_{s}$. Then $\mathrm{Var}^{S,X}_{s}(x) = \delta^{s,x}$, if $P = x$, for a unique $x\in X_{s}$; $\mathrm{Var}^{S,X}_{s}(\sigma^{\mathbf{T}_{\Sigma}(X)}) = \varnothing^{S}$, if $P = \sigma^{\mathbf{T}_{\Sigma}(X)}$, for a unique $\sigma \in \Sigma_{\lambda,s}$; and, finally, $\mathrm{Var}^{S,X}_{s}(\sigma^{\mathbf{T}_{\Sigma}(X)}((P_{j})_{j\in\bb{\mathbf{s}}})) = \bigcup_{j\in\bb{\mathbf{s}}}\mathrm{Var}^{S,X}_{s_{j}}(P_{j})$, if $P = \sigma^{\mathbf{T}_{\Sigma}(X)}((P_{j})_{j\in\bb{s}})$, for a unique $\mathbf{s}\in S^{\star}-\{\lambda\}$, a unique $\sigma \in \Sigma_{\mathbf{s},s}$ and a unique family $(P_{j})_{j\in\bb{\mathbf{s}}}\in \mathrm{T}_{\Sigma}(X)_{\mathbf{s}}$.
\end{remark}

\begin{figure}
$$
\xymatrix{ X \ar[r]^-{\eta^{X}}
  \ar[rd]_-{\delta^{S,X}}
 &
\mathrm{T}_{\Sigma}(X)
  \ar[d]^{\ext{(\delta^{S,X})} = \mathrm{Var}^{S,X} = (\mathrm{Var}^{S,X}_{s})_{s\in S}} \\
 &
{\mathrm{Fin}^{S}}(X) = (\mathrm{Sub}_{\mathrm{f}}(X))_{s\in S}}
$$
\caption{Variables.}
\label{FVar}
\end{figure}

\section{Derived operations}

Next we define, for an $S$-sorted signature $\Sigma$, a $\Sigma$-algebra $\mathbf{A}$ and a word $\mathbf{s}$ on $S$, the concept of many-sorted $\mathbf{s}$-ary operation on $\mathbf{A}$, that of many-sorted $\mathbf{s}$-ary derived operation on $\mathbf{A}$, and the procedure of realization of terms $P$ of type $(\vs{\mathbf{s}},s)$ as derived operations $P^{\mathbf{A}}$ on $\mathbf{A}$ (for the meaning of the phrase  ``term of type $(\vs{\mathbf{s}},s)$'' see below). But before doing that, we investigate the relationships between some functors defined over a category associated with the free monoid on the set of sorts $S$.

\begin{definition}
Let $\mathbf{S}^{\star}$ be the free monoid on $S$ and $\mathsf{C}(\mathbf{S}^{\star})$ the category whose objects are the words $\mathbf{s}$ on $S$ and whose morphisms from $\mathbf{s}$ to $\mathbf{s}'$ are the isotone and injective mappings 
$\varphi$ from $\bb{\mathbf{s}}$ to $\bb{\mathbf{s}'}$ such that $\mathbf{s}'\circ \varphi= \mathbf{s}$. Then we let $L$ stand for the contravariant functor from  
$\mathsf{C}(\mathbf{S}^{\star})$ to $\mathsf{Set}^{\mathsf{Set}^{S}}$ whose object mapping sends $\mathbf{s}$ to the functor $L(\mathbf{s})$ from $\mathsf{Set}^{S}$ to 
$\mathsf{Set}$ that sends an $S$-sorted set $A$ to $A_{\mathbf{s}} = \prod_{i\in \bb{\mathbf{s}}}A_{s_{i}}$ and an $S$-sorted mapping $f\colon A\mor B$ to the unique mapping $\prod_{i\in \bb{\mathbf{s}}}f_{s_{i}}$, denoted simply as $f_{\mathbf{s}}$, from $A_{\mathbf{s}}$ to $B_{\mathbf{s}}$ such that, for every $i\in \bb{\mathbf{s}}$, the following diagram
\begin{center}
\begin{tikzpicture}
[ACliment/.style={-{To [angle'=45, length=5.75pt, width=4pt, round]}},
RHACliment/.style={right hook-{To [angle'=45, length=5.75pt, width=4pt, round]}, font=\scriptsize}, 
scale=1
]

\node[] (AS) at (0,0) {$A_{\mathbf{s}}$};
\node[] (As) at (0,-2) {$A_{s_{i}}$};
\node[] (BS) at (4,0) {$B_{\mathbf{s}}$};
\node[] (Bs) at (4,-2) {$B_{s_{i}}$};

\draw[ACliment]  (AS) to node [above ] {$f_{\mathbf{s}}$} (BS);
\draw[ACliment]  (As) to node [below ] {$f_{s_{i}}$} (Bs);
\draw[ACliment]  (AS) to node [left ] {$\mathrm{pr}_{i}$} (As);
\draw[ACliment]  (BS) to node [right ] {$\mathrm{pr}_{i}$} (Bs);
\end{tikzpicture}
\end{center}
commutes; and whose morphism mapping sends $\varphi\colon \mathbf{s}\mor \mathbf{s}'$ to the natural transformation $L(\varphi)$ from $L(\mathbf{s}')$ to $L(\mathbf{s})$ that assigns to an $S$-sorted set $A$ the mapping $L(\varphi)_{A}$ from $A_{\mathbf{s}'}$ to $A_{\mathbf{s}}$ that sends $a\in A_{\mathbf{s}'}$ to $a\circ\varphi = a_{\varphi}\in A_{\mathbf{s}}$. 

Moreover, we let $R$ stand for the contravariant functor from $\mathsf{C}(\mathbf{S}^{\star})$ to $\mathsf{Set}^{\mathsf{Set}^{S}}$ whose object mapping sends $\mathbf{s}$ to the functor $R(\mathbf{s})$ from $\mathsf{Set}^{S}$ to $\mathsf{Set}$ that sends an $S$-sorted set $A$ to $A_{\vs{\mathbf{s}}} = \mathrm{Hom}(\vs{\mathbf{s}},A)$, where 
$\vs{\mathbf{s}}$ is the $S$-sorted set defined, for every $s\in S$, as 
$(\vs{\mathbf{s}})_{s} = \mathbf{s}^{-1}[\{s\}]$, and an $S$-sorted mapping $f\colon A\mor B$ to the mapping 
$\mathrm{Hom}(\vs{\mathbf{s}},f)$, also denoted $f_{\vs{\mathbf{s}}}$, from $\mathrm{Hom}(\vs{\mathbf{s}},A)$ to $\mathrm{Hom}(\vs{\mathbf{s}},B)$ that sends $a\in \mathrm{Hom}(\vs{\mathbf{s}},A)$ to $f\circ a\in \mathrm{Hom}(\vs{\mathbf{s}},B)$, i.e., $R(\mathbf{s})$ is the hom functor $\mathrm{Hom}(\vs{\mathbf{s}},\bigcdot)$; and whose morphism mapping sends $\varphi\colon \mathbf{s}\mor \mathbf{s}'$ to the natural transformation $R(\varphi)$ from $R(\mathbf{s}')$ to $R(\mathbf{s})$ that assigns to an $S$-sorted set $A$ the mapping $R(\varphi)_{A}$ from $A_{\vs{\mathbf{s}'}}$ to $A_{\vs{\mathbf{s}}}$ that sends $a\in A_{\vs{\mathbf{s}'}}$ to $a\circ\varphi = a_{\varphi}\in A_{\vs{\mathbf{s}}}$.
\end{definition}

\begin{remark}
If $\varphi$ is a morphism from $\mathbf{s}$ to $\mathbf{s}'$, then $\mathbf{s}'$ can be represented as:
$$ 
\mathbf{s}' = \mathbf{u}_{0}\bconcat s_{0}\bconcat \mathbf{u}_{1}
\bconcat\ldots \bconcat
\mathbf{u}_{\bb{\mathbf{s}}-1}\bconcat s_{\bb{\mathbf{s}}-1}\bconcat \mathbf{u}_{\bb{\mathbf{s}}},
$$
for a univocally determined family $(\mathbf{u}_{i})_{i\in \bb{\mathbf{s}}+1}$ of words on $S$ (and where, for every $i\in \bb{\mathbf{s}}$, $s_{i} = s'_{\varphi(i)}$). And, reciprocally, if for two words $\mathbf{s}$ and 
$\mathbf{s}'$, the word $\mathbf{s}'$ can be represented as a concatenation of words as just said and satisfying the letters $s_{i}$ appropriate conditions, then there exists a unique morphism from $\mathbf{s}$ to $\mathbf{s}'$ which induces the given representation of $\mathbf{s}'$.
\end{remark}

\begin{proposition}\label{isovs}
Let $\mathbf{s}$ be an element of $S^{\star}$ and $A$ an $S$-sorted set. Then there exists an isomorphism $\theta_{\mathbf{s},A}$ from $A_{\vs{\mathbf{s}}}$
to $A_{\mathbf{s}}$. Moreover, for every $S$-sorted mapping $f\colon
A\mor B$, the following diagram
\begin{center}
\begin{tikzpicture}
[ACliment/.style={-{To [angle'=45, length=5.75pt, width=4pt, round]}},
RHACliment/.style={right hook-{To [angle'=45, length=5.75pt, width=4pt, round]}, font=\scriptsize}, 
scale=1
]

\node[] (RA) at (0,0) {$A_{\vs{\mathbf{s}}}$};
\node[] (RB) at (0,-2) {$B_{\vs{\mathbf{s}}}$};
\node[] (LA) at (4,0) {$A_{\mathbf{s}}$};
\node[] (LB) at (4,-2) {$B_{\mathbf{s}}$};

\draw[ACliment]  (RA) to node [above ] {$\theta_{\mathbf{s},A}$} (LA);
\draw[ACliment]  (RB) to node [below ] {$\theta_{\mathbf{s},B}$} (LB);
\draw[ACliment]  (RA) to node [left ] {$f_{\vs{\mathbf{s}}}$} (RB);
\draw[ACliment]  (LA) to node [right ] {$f_{\mathbf{s}}$} (LB);
\end{tikzpicture}
\end{center}
commutes. Therefore the family $\theta = (\theta_{\mathbf{s}})_{\mathbf{s}\in S^{\star}}$ is a natural isomorphism from the functor $R$ to the functor $L$. 
\end{proposition}


\begin{definition}
We let $\vs{\,\bigcdot\,}$ stand for the functor from $\mathsf{C}(\mathbf{S}^{\star})$ to $\mathsf{Set}^{S}$ whose object mapping sends $\mathbf{s}$ to $\vs{\mathbf{s}}$ and whose morphism mapping sends $\varphi\colon \mathbf{s}\mor \mathbf{s}'$ to the $S$-sorted mapping $\vs{\varphi}$ from $\vs{\mathbf{s}}$ to $\vs{\mathbf{s}'}$ defined, for every $s\in S$ and every $i\in (\vs{\mathbf{s}})_{s}$, as  $(\vs{\varphi})_{s}(i) = \varphi(i)$. 

Let $V^{S} = (V^{S}_{s})_{s\in S}$ be an arbitrary, but fixed, $S$-sorted set of variables with a countably infinity of variables in each component, say    $V^{S}_{s} = \{\,v^{s}_{n}\mid n\in \mathbb{N}\,\}$, for every $s\in S$. Then we let $V^{S}_{\,\bigcdot\,}$ stand for the functor from $\mathsf{C}(\mathbf{S}^{\star})$ to $\mathsf{Set}^{S}$ whose object mapping sends $\mathbf{s}$ to the finite subset $V^{S}_{\mathbf{s}}$ of $V^{S}$, where $V^{S}_{\mathbf{s}}$ is defined, for every $s\in S$, as
$
V^{S}_{\mathbf{s},s} = \{\,v^{s}_{i}\in V^{S}_{s}\mid i\in \bb{\mathbf{s}}_{s}\,\}
$ 
(observe that $V^{S}_{\mathbf{s},s}$ is empty for those sorts $s$ that do
not occur in the word $\mathbf{s}$); and whose morphism mapping sends $\varphi\colon \mathbf{s}\mor \mathbf{s}'$ to the $S$-sorted mapping $V^{S}_{\varphi}$ from 
$V^{S}_{\mathbf{s}}$ to $V^{S}_{\mathbf{s}'}$ defined, for every $s\in S$ and every 
$j\in \bb{\mathbf{s}}_{s}$, as $V^{S}_{\varphi,s}(v^{s}_{j}) = 
v^{s}_{j'}$, where, by Definition~\ref{DSubw}, for  $(i_{j})_{j\in\bb{\mathbf{s}}_{s}}$, the enumeration in ascending order of the occurrences of $s$ in $\mathbf{s}$, and $(i_{k})_{k\in\bb{\mathbf{s}'}_{s}}$, the enumeration in ascending order of the occurrences of $s$ in $\mathbf{s}'$, $j'$ is the unique element in $\bb{\mathbf{s}'}_{s}$ such that  $\varphi(i_{j}) = i_{j'}$.
\end{definition}

\begin{proposition}\label{isovsV}
Let $\mathbf{s}$ be an element of $S^{\star}$. Then there exists an isomorphism 
$\vartheta_{\mathbf{s}}$ from $\vs{\mathbf{s}}$
to $V^{S}_{\mathbf{s}}$. Moreover, for every morphism $\varphi\colon
\mathbf{s}\mor \mathbf{s}'$, the following diagram
\begin{center}
\begin{tikzpicture}
[ACliment/.style={-{To [angle'=45, length=5.75pt, width=4pt, round]}},
RHACliment/.style={right hook-{To [angle'=45, length=5.75pt, width=4pt, round]}, font=\scriptsize}, 
scale=1
]

\node[] (s) at (0,0) {$\vs{\mathbf{s}}$};
\node[] (s') at (0,-2) {$\vs{\mathbf{s}'}$};
\node[] (Vs) at (4,0) {$V^{S}_{\mathbf{s}}$};
\node[] (Vs') at (4,-2) {$V^{S}_{\mathbf{s}'}$};

\draw[ACliment]  (s) to node [above ] {$\vartheta_{\mathbf{s}}$} (Vs);
\draw[ACliment]  (s') to node [below ] {$\vartheta_{\mathbf{s}'}$} (Vs');
\draw[ACliment]  (s) to node [left ] {$\vs{\varphi}$} (s');
\draw[ACliment]  (Vs) to node [right ] {$V^{S}_{\varphi}$} (Vs');
\end{tikzpicture}
\end{center}
commutes. Therefore the family $\vartheta = (\vartheta_{\mathbf{s}})_{\mathbf{s}\in S^{\star}}$ is a natural isomorphism from the functor $\vs{\,\bigcdot\,}$ to the functor $V^{S}_{\,\bigcdot\,}$. 
\end{proposition}

\begin{remark}
In virtue of the previous propositions, for every word $\mathbf{s}\in S^{\star}$ and every $S$-sorted set $A$, the sets $A_{\vs{\mathbf{s}}}$, $A_{\mathbf{s}}$ and $A_{V^{S}_{\mathbf{s}}}$ are naturally  isomorphic.
In what follows, if there is no ambiguity, we will not distinguish notationally between the elements of $A_{\vs{\mathbf{s}}}$ and those of $A_{\mathbf{s}}$ or those of $A_{V^{S}_{\mathbf{s}}}$. Also, for simplicity of notation, sometimes, instead of using $V^{S}_{\mathbf{s}}$ we will use 
$\vs{\mathbf{s}}$.
\end{remark}

\begin{definition}\label{wdop}
Let $\mathbf{s}$ be an element of $S^{\star}$, $\mathbf{A}$ a $\Sigma$-algebra, $s\in S$ and $P\in \mathrm{T}_{\Sigma}(\vs{\mathbf{s}})_{s}$ a term for $\Sigma$ of type $(\vs{\mathbf{s}},s)$. Then
\begin{enumerate}
\item the $\Sigma$-algebra of the \emph{many-sorted} $\vs{\mathbf{s}}$-\emph{ary} \emph{operations on} $\mathbf{A}$, $\mathbf{O}_{\vs{\mathbf{s}}}(\mathbf{A})$, is 
$\mathbf{A}^{A_{\vs{\mathbf{s}}}}$, i.e., the direct $A_{\vs{\mathbf{s}}}$-power of 
$\mathbf{A}$, where $A_{\vs{\mathbf{s}}}$ is $\mathrm{Hom}(\vs{\mathbf{s}},A)$, the (ordinary) set of the $S$-sorted mappings from $\vs{\mathbf{s}}$ to $A$ (which is naturally isomorphic to $A_{\mathbf{s}}$). From now on, to shorten terminology, we will speak of $\mathbf{s}$-\emph{ary} \emph{operations on} $\mathbf{A}$ instead of \emph{many-sorted} $\vs{\mathbf{s}}$-\emph{ary} \emph{operations on} $\mathbf{A}$, and, for brevity, we will write $\mathbf{O}_{\mathbf{s}}(\mathbf{A})$ instead of $\mathbf{O}_{\vs{\mathbf{s}}}(\mathbf{A})$;

\item the $\Sigma$-algebra of the \emph{many-sorted} $\mathbf{s}$-\emph{ary} \emph{derived operations on} $\mathbf{A}$, $\mathbf{D}_{\vs{\mathbf{s}}}(\mathbf{A})$, is the subalgebra of $\mathbf{O}_{\vs{\mathbf{s}}}(\mathbf{A})$ generated by the subfamily
$$
\mathcal{P}^{A}_{\vs{\mathbf{s}}} = (\mathcal{P}^{A}_{\vs{\mathbf{s}},s})_{s\in S} =
(\{\,\mathrm{pr}^{A}_{\vs{\mathbf{s}},s,i} \mid i\in (\vs{\mathbf{s}})_{s}\,\})_{s\in S}
$$
of $\mathrm{O}_{\vs{\mathbf{s}}}(A) = A^{A_{\vs{\mathbf{s}}}}$, where, for every $s\in S$ and $i\in (\vs{\mathbf{s}})_{s}$, $\mathrm{pr}^{A}_{\vs{\mathbf{s}},s,i}$ is the mapping from $A_{\vs{\mathbf{s}}}$ to $A_{s}$ which sends $a\in A_{\vs{\mathbf{s}}}$ to $a_{s}(i)$ (note that the set $\{\,\mathrm{pr}^{A}_{\vs{\mathbf{s}},s,i} \mid i\in (\vs{\mathbf{s}})_{s}\,\}$ has $\mathrm{card}((\vs{\mathbf{s}})_{s})$ projections if $A_{\vs{\mathbf{s}}} \neq \varnothing$ and 
$\mathrm{card}(A_{s}) > 1$, or a unique mapping if $A_{\vs{\mathbf{s}}} = \varnothing$ or $\mathrm{card}(A_{s}) = 1$). From now on, to shorten terminology, we will speak of $\mathbf{s}$-\emph{ary} \emph{derived operations on} $\mathbf{A}$ instead of \emph{many-sorted} 
$\vs{\mathbf{s}}$-\emph{ary} \emph{derived operations on} $\mathbf{A}$, and, for brevity, we will write $\mathbf{D}_{\mathbf{s}}(\mathbf{A})$ instead of $\mathbf{D}_{\vs{\mathbf{s}}}(\mathbf{A})$;

\item we denote by $\mathrm{d}^{\mathbf{s},\mathbf{A}}$ the unique homomorphism from 
$\mathbf{T}_{\Sigma}(\vs{\mathbf{s}})$ to $\mathbf{O}_{\mathbf{s}}(\mathbf{A})$ such that $\mathrm{pr}^{A}_{\vs{\mathbf{s}}} = \mathrm{d}^{\mathbf{s},\mathbf{A}}\circ \eta^{\vs{\mathbf{s}}}$, where $\mathrm{pr}^{A}_{\vs{\mathbf{s}}}$ is the $S$-sorted mapping $(\mathrm{pr}^{A}_{\vs{\mathbf{s}},s})_{s\in S}$ from $\vs{\mathbf{s}}$ to 
$\mathrm{O}_{\mathbf{s}}(A)$, with $\mathrm{pr}^{A}_{\vs{\mathbf{s}},s} = (\mathrm{pr}^{A}_{\vs{\mathbf{s}},s,i})_{i\in (\vs{\mathbf{s}})_{s}}$, for every $s\in S$. Furthermore, $P^{\mathbf{A}}$ denotes the image of $P$ under $\mathrm{d}^{\mathbf{s},\mathbf{A}}_{s}$, and we call the mapping $P^{\mathbf{A}}$ from $A_{\vs{\mathbf{s}}}$ to $A_{s}$, the \emph{derived operation on} $\mathbf{A}$ \emph{determined by} $P$, or, less frequently, the \emph{term realization of} $P$ \emph{on} $\mathbf{A}$. Moreover, to simplify the notation, we will also denote by $\mathrm{d}^{\mathbf{s},\mathbf{A}}$ the corestriction of the homomorphism $\mathrm{d}^{\mathbf{s},\mathbf{A}}\colon\mathbf{T}_{\Sigma}(\vs{\mathbf{s}})\mor \mathbf{O}_{\mathbf{s}}(\mathbf{A})$ to the subalgebra $\mathbf{D}_{\mathbf{s}}(\mathbf{A})$ of $\mathbf{O}_{\mathbf{s}}(\mathbf{A})$.
\end{enumerate}
\end{definition}

\begin{remark}
What we have called \emph{derived operations on} $\mathbf{A}$, following the terminology in Cohn~\cite{pC81}, pp. 145--149, are also known, for those following that one in Gr\"{a}tzer~\cite{gG08}, pp. 37--45, and J\'{o}nsson~\cite{bJ72}, pp. 83--87, as \emph{polynomial operations of} $\mathbf{A}$.
\end{remark}

\begin{proposition}[Reciprocity law]\label{Rlw}
Let $\mathbf{s}$ be an element of $S^{\star}$, $s\in S$, 
$P\in \mathrm{T}_{\Sigma}(\vs{\mathbf{s}})_{s}$, $\mathbf{A}$ a $\Sigma$-algebra and $a$ an $S$-sorted mapping from $\vs{\mathbf{s}}$ to $A$. Then 
$$
a_{s}^{\sharp}(P) = P^{\mathbf{A}}(a).
$$  
\end{proposition}

\section{Derivors}

In this subsection, after defining the variety of Hall algebras, we introduce the notion of derivor between many-sorted signatures. This will allow us, using the homomorphisms between Hall algebras, to obtain $\mathsf{Sig}_{\mathfrak{d}}$, the category of many-sorted signatures and derivors. Next, after defining a suitable contravariant functor from $\mathsf{Sig}_{\mathfrak{d}}$ to $\mathsf{Cat}$, the category of $\boldsymbol{\mathcal{U}}$-locally small categories and functors, we obtain, by means of the Grothendieck construction, $\mathsf{Alg}_{\mathfrak{d}}$, the category with objects the ordered pairs $((S,\Sigma),(A,F))$, with $(S,\Sigma)$ a many-sorted signature and $(A,F)$ a $\Sigma$-algebra, and morphisms from $((S,\Sigma),(A,F))$ to  $((T,\Lambda),(B,G))$ the ordered pairs $((\varphi,d),f)$, with $(\varphi,d)$ a derivor from $(S,\Sigma)$ to $(T,\Lambda)$ and $f$ a homomorphism of $\Sigma$-algebras from $(A,F)$ to $(B_{\varphi},G^{(\varphi,d)})$, a canonically derived algebra of $(B,G)$. 

Before defining the notion of Hall algebra, we recall that
\begin{enumerate}
\item a finitary specification is an ordered triple $(S,\Sigma,\mathcal{E})$, where $S$ is a set of sorts, $\Sigma$ an $S$-sorted signature, and $\mathcal{E}\subseteq\mathrm{Eq}_{V^{S}}(\Sigma)=
(\mathrm{T}_{\Sigma}(X)_{s}^{2})_{(X,s)\in \mathrm{Sub}_{\mathrm{f}}(V^{S})\times S}$, i.e., a set of finitary  $\Sigma$-equations, where $V^{S}$, the $S$-sorted set of the variables, is a fixed $S$-countably infinite $S$-sorted set; and that
\item a $\Sigma$-algebra $\mathbf{A}$ is a $(S,\Sigma,\mathcal{E})$-algebra if $\mathbf{A}\models^{\Sigma} \mathcal{E}$, i.e., if, for every $(X,s)\in \mathrm{Sub}_{\mathrm{f}}(V^{S})\times S$ and every $(P,Q)\in \mathcal{E}_{X,s}$, $\mathbf{A}$ is a model of $(P,Q)$, in symbols $\mathbf{A}\models^{\Sigma}_{X,s} (P,Q)$, which in turn means that, for every $a\in \mathrm{Hom}(X,A)$, $a^{\sharp}_{s}(P) = a^{\sharp}_{s}(Q)$.
\end{enumerate}

\begin{definition}
Let $S$ be a set of sorts and let $V^{\mathrm{H}_{S}} = (V^{\mathrm{H}_{S}}_{\mathbf{s},s})_{(\mathbf{s},s)\in S^{\star}\times S}$ be a fixed $S^{\star}\times S$-sorted set of variables with a countably infinity of variables in each component, say $V^{\mathrm{H}_{S}}_{\mathbf{s},s} = \{\,v^{\mathbf{s},s}_{n}\mid n\in \mathbb{N}\,\}$ (which is isomorphic to 
$\{(\mathbf{s},s)\}\times\mathbb{N}$), for every $(\mathbf{s},s)\in S^{\star}\times S$.  A \emph{Hall algebra for} $S$ is a $\mathrm{H}_{S} = (S^{\star}\times S,\Sigma^{\mathrm{H}_{S}},\mathcal{E}^{\mathrm{H}_{S}})$-algebra, where $\Sigma^{\mathrm{H}_{S}}$ is the $S^{\star}\times S$-sorted signature, i.e., the $(S^{\star}\times S)^{\star}\times (S^{\star}\times S)$-sorted set, defined as follows:
\begin{enumerate}
\item[$\mathrm{HS}_{1}$.] For every $\mathbf{s}\in S^{\star}$ and $i\in\bb{\mathbf{s}}$,
      $$
      \pi^{\mathbf{s}}_{i}\colon\lambda\mor(\mathbf{s},s_{i}),
      $$
      where $\bb{\mathbf{s}}$ is the \emph{length} of the word $w$ and $\lambda$ the
      \emph{empty word} in the underlying set of the free monoid on
      $S^{\star}\times S$.

\item[$\mathrm{HS}_{2}$.] For every $\mathbf{r}$, $\mathbf{s}\in S^{\star}$ and $s\in S$,
      $$
      \xi_{\mathbf{r},\mathbf{s},s}\colon ((\mathbf{s},s),(\mathbf{r},s_{0}),\ldots,(\mathbf{r},s_{\bb{\mathbf{s}}-1}))\mor(\mathbf{r},s);
      $$
\end{enumerate}
while $\mathcal{E}^{\mathrm{H}_{S}}$ is the sub-$(S^{\star}\times
S)^{\star}\times (S^{\star}\times S)$-sorted set of
$\mathrm{Eq}(\Sigma^{\mathrm{H}_{S}})$, where
$$
\mathrm{Eq}(\Sigma^{\mathrm{H}_{S}}) =
(\mathrm{T}_{\Sigma^{\mathrm{H}_{S}}}(\vs{\ol{\mathbf{s}}})_{(\mathbf{r},s)}^{2})_{(\ol{\mathbf{s}},(\mathbf{r},s))\in
(S^{\star}\times S)^{\star}\times (S^{\star}\times S)},
$$
defined as follows:
\begin{enumerate}
\item[$\mathrm{H}_{1}$.] \emph{Projection}.
      For every $\mathbf{r}$, $\mathbf{s}\in S^{\star}$ and $i\in\bb{\mathbf{s}}$, the equation
      $$
      \xi_{\mathbf{r},\mathbf{s},s_{i}}(\pi^{\mathbf{s}}_{i},v^{\mathbf{r},s_{0}}_{0},\ldots,
           v^{\mathbf{r},s_{\bb{\mathbf{s}}-1}}_{\bb{\mathbf{s}}-1})=v^{\mathbf{r},s_{i}}_{i}
      $$
      of type
      $(((\mathbf{r},s_{0}),\ldots,(\mathbf{r},s_{\bb{\mathbf{s}}-1})),(\mathbf{r},s_{i})).$
\item[$\mathrm{H}_{2}$.] \emph{Identity}.
      For every $\mathbf{r}\in S^{\star}$ and $j\in \bb{\mathbf{r}}$, the equation
      $$
      \xi_{\mathbf{r},\mathbf{r},r_{j}}(v^{\mathbf{r},r_{j}}_{j},\pi^{\mathbf{r}}_{0},\ldots,\pi^{\mathbf{r}}_{\bb{\mathbf{r}}-1})=
           v^{\mathbf{r},r_{j}}_{j}
      $$
      of type
      $(((\mathbf{r},r_{j})),(\mathbf{r},r_{j})).$
\item[$\mathrm{H}_{3}$.] \emph{Associativity}.
      For every $\mathbf{q}$, $\mathbf{r}$, $\mathbf{s}\in S^{\star}$ and $s\in S$, the
      equation
\begin{multline*}
\xi_{\mathbf{q},\mathbf{r},s}(
      \xi_{\mathbf{r},\mathbf{s},s}(v^{\mathbf{s},s}_{0},v^{\mathbf{r},s_{0}}_{1},\ldots,
          v^{\mathbf{r},s_{\bb{\mathbf{s}}-1}}_{\bb{\mathbf{s}}}),
          v^{\mathbf{q},r_{0}}_{\bb{\mathbf{s}}+1},
          \ldots,v^{\mathbf{q},r_{\bb{\mathbf{r}}-1}}_{\bb{\mathbf{s}}+\bb{\mathbf{r}}}) = \\
       \xi_{\mathbf{q},\mathbf{s},s}(v^{\mathbf{s},s}_{0},
          \xi_{\mathbf{q},\mathbf{r},s_{0}}(v^{\mathbf{r},s_{0}}_{1},v^{\mathbf{q},r_{0}}_{\bb{\mathbf{s}}+1},
                   \ldots,v^{\mathbf{q},r_{\bb{\mathbf{r}}-1}}_{\bb{\mathbf{s}}+\bb{\mathbf{r}}}),
                   \\
          \ldots, 
       \xi_{\mathbf{q},\mathbf{r},s_{\bb{\mathbf{s}}-1}}(v^{\mathbf{r},s_{\bb{\mathbf{s}}-1}}_{\bb{\mathbf{s}}},
                   v^{\mathbf{q},r_{0}}_{\bb{\mathbf{s}}+1},
       \ldots,v^{\mathbf{q},r_{\bb{\mathbf{r}}-1}}_{\bb{\mathbf{s}}+\bb{\mathbf{r}}}))
\end{multline*}
     of type
     $(((\mathbf{s},s),(\mathbf{r},s_{0}),\ldots,(\mathbf{r},s_{\bb{\mathbf{s}}-1}),
            (\mathbf{q},r_{0}),\ldots,(\mathbf{q},r_{\bb{\mathbf{r}}-1})),(\mathbf{q},s)).$
\end{enumerate}

We call the formal constants $\pi^{\mathbf{s}}_{i}$ \emph{projections}, and the formal operations $\xi_{\mathbf{r},\mathbf{s},s}$ \emph{substitution operators}. Furthermore, we will denote by $\mathsf{Alg}(\mathrm{H}_{S})$ the category of Hall algebras for $S$ and homomorphisms between Hall algebras.  Since $\mathsf{Alg}(\mathrm{H}_{S})$ is a variety, the forgetful functor $\mathrm{G}_{\mathrm{H}_{S}}$ from $\mathsf{Alg}(\mathrm{H}_{S})$ to $\mathsf{Set}^{S^{\star}\times S}$ has a left adjoint $\mathbf{T}_{\mathrm{H}_{S}}$, situation denoted by $\mathbf{T}_{\mathrm{H}_{S}}\dashv \mathrm{G}_{\mathrm{H}_{S}}$, or diagrammatically by
$$
\xymatrix@=5pc{ \mathsf{Alg}(\mathrm{H}_{S})
\ar@<1.5ex>[r]^{\mathrm{G}_{\mathrm{H}_{S}}} \ar@{}[r]|{\uadj} &
\mathsf{Set}^{S^{\star}\times S}
\ar@<1.5ex>[l]^{\mathbf{T}_{\mathrm{H}_{S}}} }
$$
which assigns to an $S^{\star}\times S$-sorted set $\Sigma$ the corresponding free Hall algebra
$\mathbf{T}_{\mathrm{H}_{S}}(\Sigma)$.
\end{definition}

\begin{remark}
The category $\mathsf{Alg}(\mathrm{H}_{S})$, of Hall algebras for $S$, is an algebraic version of many-sorted clone theory.
\end{remark}

\begin{remark}
From $\mathrm{H}_{3}$, for $\mathbf{s}=\lambda$, the empty word on $S$, we get the invariance of constant functions axiom in~\cite{gm85}: For every $\mathbf{q}$, $\mathbf{r}\in S^{\star}$ and $s\in S$, we have the equation
$$
\xi_{\mathbf{q},\mathbf{r},s}(\xi_{\mathbf{r},\lambda,s}(v^{\lambda,s}_{0}),v^{\mathbf{q},r_{0}}_{1},
\ldots,v^{\mathbf{q},r_{\bb{\mathbf{r}}-1}}_{\bb{\mathbf{r}}})=
\xi_{\mathbf{q},\lambda,s}(v^{\lambda,s}_{0})
$$
of type $(((\lambda,s),(\mathbf{q},r_{0}),\ldots,(\mathbf{q},r_{\bb{\mathbf{r}}-1})),(\mathbf{q},s))$.
\end{remark}

For every $S$-sorted set $A$, $\mathrm{O}_{\mathrm{H}_{S}}(A) = (\mathrm{Hom}(A_{\mathbf{s}},A_{s}))_{(\mathbf{s},s)\in S^{\star}\times S}$, the $S^{\star}\times S$-sorted set of operation for $A$, is naturally equipped with a structure of Hall algebra, as stated in the following proposition, if we realize the projections as the true projections and the substitution operators as the generalized composition of mappings.

\begin{proposition}
Let $A$ be an $S$-sorted set and $\mathbf{O}_{\mathrm{H}_{S}}(A)$ the $\Sigma^{\mathrm{H}_{S}}$-algebra with underlying many-sorted set $\mathrm{O}_{\mathrm{H}_{S}}(A)$ and algebraic structure defined as follows:
\begin{enumerate}
\item For every $\mathbf{s}\in S^{\star}$ and $i\in\bb{\mathbf{s}}$,
      $(\pi^{\mathbf{s}}_{i})^{\mathbf{O}_{\mathrm{H}_{S}}(A)} =
      \mathrm{pr}^{A}_{\mathbf{s},i}\colon A_{\mathbf{s}}\mor A_{s_{i}}$.

\item For every $\mathbf{r},\mathbf{s}\in S^{\star}$ and $s\in S$,
      $\xi_{\mathbf{r},\mathbf{s},s}^{\mathbf{O}_{\mathrm{H}_{S}}(A)}$ is defined, for every $f\in
      A_{s}^{ A_{\mathbf{s}} }$ and $g\in A^{A_{\mathbf{r}}}_{\mathbf{s}}$, as
      $\xi_{\mathbf{r},\mathbf{s},s}^{\mathbf{O}_{\mathrm{H}_{S}}(A)}(f,g_{0},\ldots,g_{\bb{\mathbf{s}}-1})
      = f\circ\langle g_{i}\rangle_{i\in\bb{\mathbf{s}} }$, where $\langle
      g_{i}\rangle_{i\in\bb{\mathbf{s}}}$ is the unique mapping from $A_{\mathbf{r}}$ to
      $A_{\mathbf{s}}$ such that, for every $i\in \bb{\mathbf{s}}$, we have that
      $$
      \mathrm{pr}^{A}_{\mathbf{s},i}\circ\langle g_{i}\rangle_{i\in\bb{\mathbf{s}}} = g_{i}.
      $$
\end{enumerate}
Then $\mathbf{O}_{\mathrm{H}_{S}}(A)$ is a Hall algebra, the \emph{Hall algebra for} $(S,A)$.
\end{proposition}

\begin{remark}
The closed sets of the Hall algebra $\mathbf{O}_{\mathrm{H}_{S}}(A)$ for $(S,A)$ are precisely the clones of (many-sorted) operations on the $S$-sorted set $A$. On the other hand, every $\Sigma$-algebra $\mathbf{A}$ has associated a Hall algebra. In fact, it suffices to consider $\mathbf{O}_{\mathrm{H}_{S}}(A)$, denoted by $\mathbf{O}_{\mathrm{H}_{S}}(\mathbf{A})$. Moreover, the finitary derived operations on $\mathbf{A}$ is a subalgebra of the Hall algebra $\mathbf{O}_{\mathrm{H}_{S}}(\mathbf{A})$.
\end{remark}

For every $S$-sorted signature $\Sigma$, $\mathrm{Tm}_{\mathrm{H}_{S}}(\Sigma) = (\mathrm{T}_{\Sigma}(\vs{\mathbf{s}})_{s})_{(\mathbf{s},s)\in S^{\star}\times S}$ is also e\-quipped with a structure of Hall algebra that formalizes the concept of substitution as stated in the following proposition.

\begin{proposition}
Let $\Sigma$ be an $S$-sorted signature and $\mathbf{Tm}_{\mathrm{H}_{S}}(\Sigma)$ the $\Sigma^{\mathrm{H}_{S}}$-algebra
with underlying many-sorted set $\mathrm{Tm}_{\mathrm{H}_{S}}(\Sigma)$ and algebraic structure defined as follows:
\begin{enumerate}
\item For every $\mathbf{s}\in S^{\star}$ and $i\in\bb{\mathbf{s}}$,
      $(\pi^{\mathbf{s}}_{i})^{\mathbf{Tm}_{\mathrm{H}_{S}}(\Sigma)}$
      is the image under $\eta^{\vs{\mathbf{s}}}_{s_{i}}$ of
      the variable $v_{i}^{s_{i}}$, where $\eta^{\vs{\mathbf{s}}} =
      (\eta^{\vs{\mathbf{s}}}_{s})_{s\in S}$ is the canonical embedding
      of $\vs{\mathbf{s}}$ into $\mathrm{T}_{\Sigma}(\vs{\mathbf{s}})$. Sometimes,
      to abbreviate, we will write $\pi^{\mathbf{s}}_{i}$ instead of
      $(\pi^{\mathbf{s}}_{i})^{\mathbf{Tm}_{\mathrm{H}_{S}}(\Sigma)}$.

\item For every $\mathbf{r},\mathbf{s}\in S^{\star}$ and $s\in S$,
      $\xi_{\mathbf{r},\mathbf{s},s}^{\mathbf{Tm}_{\mathrm{H}_{S}}(\Sigma)}$
      is the mapping
      $$\xi_{\mathbf{r},\mathbf{s},s}^{\mathbf{Tm}_{\mathrm{H}_{S}}(\Sigma)}\nfunction
      {\mathrm{T}_{\Sigma}(\vs{\mathbf{s}})_{s} \times
      \mathrm{T}_{\Sigma}(\vs{\mathbf{r}})_{s_{0}} \times \cdots \times
      \mathrm{T}_{\Sigma}(\vs{\mathbf{r}})_{s_{\bb{\mathbf{s}}-1}}}
      {\mathrm{T}_{\Sigma}(\vs{\mathbf{r}})_{s}}
      {(P,(Q_{i})_{i\in\bb{\mathbf{s}}})}
      {\left(\left(\!\begin{smallmatrix}v^{s_{i}}_{i}\\Q_{i}\end{smallmatrix}\!\right)_{i\in \bb{\mathbf{s}}}\right)^{\sharp}_{s}(P)}
     $$
     where, for $\left(\!\begin{smallmatrix}v^{s_{i}}_{i}\\Q_{i}\end{smallmatrix}\!\right)_{i\in \bb{\mathbf{s}}}$, the $S$-sorted mapping from $\vs{\mathbf{s}}$
     to $\mathrm{T}_{\Sigma}(\vs{\mathbf{r}})$ canonically associated to the family
     $(Q_{i})_{i\in\bb{\mathbf{s}}}$, $\left(\left(\!\begin{smallmatrix}v^{s_{i}}_{i}\\Q_{i}\end{smallmatrix}\!\right)_{i\in \bb{\mathbf{s}}}\right)^{\sharp}$, also denoted by $\mathcal{Q}^{\sharp}$, is the
     unique homomorphism from $\mathbf{T}_{\Sigma}(\vs{\mathbf{s}})$ into
     $\mathbf{T}_{\Sigma}(\vs{\mathbf{r}})$ such that
     $$
     \left(\left(\!\begin{smallmatrix}v^{s_{i}}_{i}\\Q_{i}\end{smallmatrix}\!\right)_{i\in \bb{\mathbf{s}}}\right)^{\sharp}\circ \eta^{\vs{\mathbf{s}}} = \left(\!\begin{smallmatrix}v^{s_{i}}_{i}\\Q_{i}\end{smallmatrix}\!\right)_{i\in \bb{\mathbf{s}}}.
     $$
     Sometimes, to abbreviate, we will write $\xi_{\mathbf{r},\mathbf{s},s}$ instead of
     $\xi_{\mathbf{r},\mathbf{s},s}^{\mathbf{Tm}_{\mathrm{H}_{S}}(\Sigma)}$.
\end{enumerate}

Then $\mathbf{Tm}_{\mathrm{H}_{S}}(\Sigma)$ is a Hall algebra, the \emph{Hall algebra for} $(S,\Sigma)$.
\end{proposition}

\begin{remark}
For every $\Sigma$-algebra $\mathbf{A}$ there exists a homomorphism of Hall algebras $\mathrm{d}^{\mathbf{A}} = (\mathrm{d}^{\mathbf{A}}_{\mathbf{\mathbf{s}},s})_{(\mathbf{\mathbf{s}},s)\in S^{\star}\times S}$ from the Hall algebra $\mathbf{Tm}_{\mathrm{H}_{S}}(\Sigma)$ to the Hall algebra $\mathbf{O}_{\mathrm{H}_{S}}(\mathbf{A})$ and its image is $\mathbf{D}_{\mathrm{H}_{S}}(\mathbf{A})$, the Hall subalgebra of $\mathbf{O}_{\mathrm{H}_{S}}(\mathbf{A})$ of the finitary derived operations on $\mathbf{A}$, which, for every $(\mathbf{\mathbf{s}},s)\in S^{\star}\times S$, is such that 
$\mathrm{D}_{\mathrm{H}_{S}}(\mathbf{A})_{\mathbf{\mathbf{s}},s} = 
\mathrm{D}_{\mathbf{\mathbf{s}}}(\mathbf{A})_{s}$, where $\mathrm{D}_{\mathbf{\mathbf{s}}}(\mathbf{A})_{s}$ is the $s$-th coordinate of the underlying $S$-sorted of 
$\mathbf{D}_{\mathbf{\mathbf{s}}}(\mathbf{A})$, the $\Sigma$-algebra of the $\mathbf{\mathbf{s}}$-\emph{ary} derived operations on $\mathbf{A}$, introduced in Definition~\ref{wdop}.
\end{remark}



Our next goal is to prove that, for every $S^{\star}\times S$-sorted set $\Sigma$, $\mathbf{T}_{\mathrm{H}_{S}}(\Sigma)$, the free Hall algebra on $\Sigma$, is isomorphic to $\mathbf{Tm}_{\mathrm{H}_{S}}(\Sigma)$. We remark that the existence of this isomorphism is interesting because it enables us to get a more tractable description of the terms in $\mathbf{T}_{\mathrm{H}_{S}}(\Sigma)$.

To attain the goal just stated we begin by defining, for a Hall algebra $\mathbf{A}$, an $S$-sorted signature $\Sigma$, an
$S^{\star}\times S$-mapping $f\colon\Sigma\mor A$, and a word $\mathbf{r}\in S^{\star}$, the concept of derived $\Sigma$-algebra of
$\mathbf{A}$ for $(f,\mathbf{r})$, since it will be used afterwards in the proof of the isomorphism between
$\mathbf{T}_{\mathrm{H}_{S}}(\Sigma)$ and $\mathbf{Tm}_{\mathrm{H}_{S}}(\Sigma)$.

\begin{definition}
Let $\mathbf{A}$ be a Hall algebra and $\Sigma$ an $S$-sorted signature. Then, for every $f\colon\Sigma\mor A$ and $\mathbf{r}\in
S^{\star}$, $\mathbf{A}^{f,\mathbf{r}}$, the \emph{derived} $\Sigma$-\emph{algebra} \emph{of} $\mathbf{A}$ \emph{for} $(f,\mathbf{r})$,
is the $\Sigma$-algebra with underlying $S$-sorted set $A^{f,\mathbf{r}} = (A_{\mathbf{r},s})_{s\in S}$ and algebraic structure $F^{f,\mathbf{r}}$, defined, for every $(\mathbf{s},s)\in S^{\star}\times S$, as
$$
F^{f,\mathbf{r}}_{\mathbf{s},s} \nfunction
{\Sigma_{\mathbf{s},s}}{\mathrm{O}_{\mathbf{s}}(A^{f,\mathbf{r}})_{s}} {\sigma}{\nfunction
    {\prod_{i\in \bb{\mathbf{s}}}A_{\mathbf{u},s_{i}}}
    {A_{\mathbf{u},s}}
    {(a_{0},\ldots,a_{\bb{\mathbf{s}}-1})}
    {\xi_{\mathbf{r},\mathbf{s},s}^{\mathbf{A}}(f_{(\mathbf{s},s)}(\sigma),a_{0},\ldots,a_{\bb{\mathbf{s}}-1})}
}
$$
where $\mathrm{O}_{\mathbf{s}}(A^{f,\mathbf{r}})_{s} = A_{\mathbf{r},s}^{\prod_{i\in \bb{\mathbf{s}}}A_{\mathbf{r},s_{i}}}$.

Furthermore, we will denote by $p^{\mathbf{r}}$ the $S$-sorted mapping from $\vs{\mathbf{r}}$ to $A^{f,\mathbf{r}}$ defined, for every $s\in S$ and $i\in\bb{\mathbf{r}}$, as $p^{\mathbf{r}}_{s}(v^{s}_{i}) = (\pi^{\mathbf{r}}_{i})^{\mathbf{A}}$, and by $(p^{\mathbf{r}})^{\sharp}$ the unique homomorphism from
$\mathbf{T}_{\Sigma}(\vs{\mathbf{r}})$ to $\mathbf{A}^{f,\mathbf{r}}$ such that $(p^{\mathbf{r}})^{\sharp}\circ \eta^{\vs{\mathbf{r}}} = p^{\mathbf{r}}$.
\end{definition}

\begin{remark}
For a $\Sigma$-algebra $\mathbf{B}=(B,G)$, we have that $G\colon\Sigma\mor\mathrm{O}_{\mathrm{H}_{S}}(B)$ and
$\mathbf{B}\iso\mathbf{O}_{\mathrm{H}_{S}}(B)^{G,\lambda}$, where $\lambda$ is the empty word on $S$.  Besides, for every $\mathbf{r}\in S^{\star}$, we have that $\mathbf{B}^{B_{\mathbf{r}}}$, the direct $B_{\mathbf{r}}$-power of $\mathbf{B}$, is isomorphic to $\mathbf{O}_{\mathrm{H}_{S}}(B)^{G,\mathbf{r}}$.
\end{remark}

\begin{lemma}\label{L:aux}
Let $\Sigma$ be an $S$-sorted signature, $\mathbf{A}$ a Hall algebra, $f\colon\Sigma\mor A$ and $\mathbf{q}\in S^{\star}$.  Then, for
every $(\mathbf{r},s)\in S^{\star}\times S$, $P\in \mathrm{T}_{\Sigma}(\vs{\mathbf{r}})_{s}$ and $a\in \prod_{i\in \bb{\mathbf{r}}}A_{\mathbf{q},r_{i}}$, we have that
$$
P^{\mathbf{A}^{f,\mathbf{q}}}(a_{0},\ldots,a_{\bb{\mathbf{r}}-1}) =
\xi_{\mathbf{q},\mathbf{r},s}^{\mathbf{A}}((p^{\mathbf{r}})^{\sharp}_{s}(P),a_{0},\ldots,a_{\bb{\mathbf{r}}-1}).
$$
\end{lemma}

\begin{proof}
By algebraic induction on the complexity of $P$.  If $P$ is a variable $v_{i}^{s}$, with $i\in\bb{\mathbf{r}}$, then
\begin{align*}
    v_{i}^{s,\mathbf{A}^{f,\mathbf{q}}}(a_{0},\ldots,a_{\bb{\mathbf{r}}-1})\ &=
    a^{\sharp}_{r_{i}}(v^{s}_{i})\\
    &=
    a_{i} \\
    &=
    \xi_{\mathbf{q},\mathbf{r},s}^{\mathbf{A}}((\pi^{\mathbf{r}}_{i})^{\mathbf{A}},a_{0},\ldots,a_{\bb{\mathbf{r}}-1} )
    \tag{by $\mathrm{H}_{1}$} \\
    &=
    \xi_{\mathbf{q},\mathbf{r},s}^{\mathbf{A}}((p^{w})^{\sharp}_{s}(v^{s}_{i}),a_{0},\ldots,a_{\bb{\mathbf{r}}-1} ).
\end{align*}
Let us assume that $P = \sigma(Q_{0},\ldots,Q_{\bb{x}-1})$, with $\sigma\in \Sigma_{\mathbf{s},s}$ and that, for every $j\in\bb{x}$, $Q_{j}\in \mathrm{T}_{\Sigma}(\vs{\mathbf{r}})_{s_{j}}$ fulfils the induction hypothesis. Then we have that
\begin{flushleft}
$(\sigma(Q_{0},\ldots,Q_{\bb{\mathbf{s}}-1}))^{\mathbf{A}^{f,\mathbf{q}}}
      (a_{0},\ldots,a_{\bb{\mathbf{r}}-1})$
\allowdisplaybreaks
\begin{align*}
    &=
    \sigma^{\mathbf{A}^{f,\mathbf{q}}}
       (
       Q_{0}^{\mathbf{A}^{f,\mathbf{q}}}(a_{0},\ldots,a_{\bb{\mathbf{r}}-1}),
       \ldots,
       Q_{\bb{\mathbf{s}}-1}^{\mathbf{A}^{f,\mathbf{q}}}(a_{0},\ldots,a_{\bb{\mathbf{r}}-1})
       )\\
    &=
    \xi_{\mathbf{q},\mathbf{s},s}^{\mathbf{A}}
       (
       f(\sigma),
       Q_{0}^{\mathbf{A}^{f,\mathbf{q}}}(a_{0},\ldots,a_{\bb{\mathbf{r}}-1}),
       \ldots,
       Q_{\bb{\mathbf{s}}-1}^{\mathbf{A}^{f,\mathbf{q}}}(a_{0},\ldots,a_{\bb{\mathbf{r}}-1})
       )\\
   &=
    \xi_{\mathbf{q},\mathbf{s},s}^{\mathbf{A}}
       (
       f(\sigma),
      \xi_{\mathbf{q},\mathbf{r},s_{0}}^{\mathbf{A}}
            ((p^{\mathbf{r}})^{\sharp}_{s_{0}}(Q_{0}),a_{0},\ldots,a_{\bb{\mathbf{r}}-1}),
       \ldots, 
\\&\qquad\qquad\qquad\qquad\qquad\qquad
\xi_{\mathbf{q},\mathbf{r},s_{\bb{\mathbf{s}}-1}}^{\mathbf{A}}
            ((p^{\mathbf{r}})^{\sharp}_{s_{\bb{\mathbf{s}}-1}}(Q_{\bb{\mathbf{s}}-1}),a_{0},\ldots,a_{\bb{\mathbf{r}}-1})
       ) \tag{by ind.}
\\&=
    \xi_{\mathbf{q},\mathbf{r},s}^{\mathbf{A}}
       (
         \xi_{\mathbf{r},\mathbf{s},s}^{\mathbf{A}}
           (f(\sigma),
            (p^{\mathbf{r}})^{\sharp}_{s_{0}}(Q_{0}),
            \ldots,
            (p^{\mathbf{r}})^{\sharp}_{s_{\bb{\mathbf{s}}-1}}(Q_{\bb{\mathbf{s}}-1})
           ),
         a_{0},
         \ldots,
         a_{\bb{\mathbf{r}}-1}
       ) \tag{by $\mathrm{H}_{3}$}   \\
    &=
    \xi_{\mathbf{q},\mathbf{r},s}^{\mathbf{A}}
       (
       \sigma^{\mathbf{A}_{\mathbf{r}}}
          (
          (p^{\mathbf{r}})^{\sharp}_{s_{0}}(Q_{0}),
          \ldots,
          (p^{\mathbf{r}})^{\sharp}_{s_{\bb{\mathbf{s}}-1}}(Q_{\bb{\mathbf{s}}-1})
          ),
       a_{0},
       \ldots,
       a_{\bb{\mathbf{r}}-1}
       ) \\
    &=
    \xi_{\mathbf{q},\mathbf{r},s}^{\mathbf{A}}
       (
       (p^{\mathbf{r}})^{\sharp}_{s}(\sigma,Q_{0},\ldots,Q_{\bb{\mathbf{s}}-1}),
       a_{0},
       \ldots,
       a_{\bb{\mathbf{r}}-1}
       ) \\
    &=
    \xi_{\mathbf{q},\mathbf{r},s}^{\mathbf{A}}
       (
       (p^{\mathbf{r}})^{\sharp}_{s}(P),
       a_{0},
       \ldots,
       a_{\bb{\mathbf{r}}-1}
       ).
\end{align*}
\end{flushleft}
This completes the proof.
\end{proof}

Next we prove that, for every $S^{\star}\times S$-sorted set $\Sigma$, the Hall algebra for $(S,\Sigma)$ is isomorphic to the free Hall algebra on $\Sigma$.

\begin{proposition}\label{iso:FrH-TerH}
Let $\Sigma$ be an $S$-sorted signature, i.e., an $S^{\star}\times S$-sorted set.  Then the Hall algebra $\mathbf{Tm}_{\mathrm{H}_{S}}(\Sigma)$ is isomorphic to $\mathbf{T}_{\mathrm{H}_{S}}(\Sigma)$, the free Hall algebra on $\Sigma$.
\end{proposition}

\begin{proof}
It is enough to prove that $\mathbf{Tm}_{\mathrm{H}_{S}}(\Sigma)$ has the universal property of $\mathbf{T}_{\mathrm{H}_{S}}(\Sigma)$. Therefore we have to specify an $S^{\star}\times S$-sorted mapping $\eta^{\Sigma}$ from $\Sigma$ to $\mathrm{Tm}_{\mathrm{H}_{S}}(\Sigma)$ such that, for every Hall algebra $\mathbf{A}$ and $S^{\star}\times S$-sorted mapping $f$ from $\Sigma$ to $A$, there is a unique homomorphism $\widehat{f}$ from $\mathbf{Tm}_{\mathrm{H}_{S}}(\Sigma)$ to $\mathbf{A}$ such that
$\widehat{f}\circ \eta^{\Sigma} = f$.  Let $\eta^{\Sigma}$ be the $S^{\star}\times S$-sorted mapping defined, for every
$(\mathbf{s},s)\in S^{\star}\times S$, as%
$$
\eta^{\Sigma}_{\mathbf{s},s}\nfunction
{\Sigma_{\mathbf{s},s}}
{\mathrm{T}_{\Sigma}(\vs{\mathbf{s}})_{s}}
{\sigma}
{\sigma\left(v^{s_{0}}_{0},\ldots,v^{s_{\bb{\mathbf{s}}-1}}_{\bb{\mathbf{s}}-1}\right)}
$$
Let $\mathbf{A}$ be a Hall algebra, $f\colon\Sigma\mor A$ an $S^{\star}\times S$-sorted mapping and $\widehat{f}$ the $S^{\star}\times S$-sorted mapping from $\mathrm{Tm}_{\mathrm{H}_{S}}(\Sigma)$ to $A$ defined, for every $(\mathbf{s},s)\in S^{\star}\times S$, as $\widehat{f}_{\mathbf{s},s} = (p^{\mathbf{s}})^{\sharp}_{s}$, where, we recall, $(p^{\mathbf{s}})^{\sharp}$ is the unique homomorphism from $\mathbf{T}_{\Sigma}(\vs{\mathbf{s}})$ to $\mathbf{A}^{f,\mathbf{s}}$ such that $(p^{\mathbf{s}})^{\sharp}\circ \eta^{\vs{\mathbf{s}}} = p^{\mathbf{s}}$. Then $\widehat{f}$ is a homomorphism of Hall algebras, because, on the one hand, for $\mathbf{s}\in S^{\star}$ and $i\in\bb{\mathbf{s}}$ we have that%
\begin{align*}
    \widehat{f}_{\mathbf{s},s_{i}}((\pi^{\mathbf{s}}_{i})^{\mathbf{Tm}_{\mathrm{H}_{S}}(\Sigma)})
    &=
    \widehat{f}_{\mathbf{s},s_{i}}(v^{s_{i}}_{i}) \\
    &=
    p^{\mathbf{s}}_{s_{i}}(v^{s_{i}}_{i}) \\
    &=
    (\pi^{\mathbf{s}}_{i})^{\mathbf{A}},
\end{align*}
and, on the other hand, for $P\in \mathrm{T}_{\Sigma}(\vs{\mathbf{s}})_{s}$ and
$\mathcal{Q} = (Q_{i})_{i\in \bb{\mathbf{s}}}\in \mathrm{T}_{\Sigma}(\vs{\mathbf{r}})_{\mathbf{s}}$ we have that %
\begin{flushleft}
$\widehat{f}_{\mathbf{r},s}(
      \xi_{\mathbf{r},\mathbf{s},s}^{\mathbf{Tm}_{\mathrm{H}_{S}}(\Sigma)}
        (P,Q_{0},\ldots,Q_{\bb{\mathbf{s}}-1})
                   )$
\begin{align*}
\quad    &=
    (p^{\mathbf{r}})^{\sharp}_{s}(
      \mathcal{Q}^{\sharp}_{s}(P)
                     ) \\
    &=
    ((p^{\mathbf{r}})^{\sharp}\circ \mathcal{Q})^{\sharp}_{s} (P)
    \tag{$(p^{\mathbf{r}})^{\sharp}\circ\mathcal{Q}^{\sharp} =
    ((p^{\mathbf{r}})^{\sharp}\circ \mathcal{Q})^{\sharp}$}
    \\
    &=
    P^{\mathbf{A}^{f,\mathbf{r}}}
      ((p^{\mathbf{r}})^{\sharp}_{s_{0}}(Q_{0}),\ldots,
        (p^{\mathbf{r}})^{\sharp}_{s_{\bb{\mathbf{s}}-1}}(Q_{\bb{\mathbf{s}}-1})
      ) \\
    &=
    \xi_{\mathbf{r},\mathbf{s},s}^{\mathbf{A}}
      ((p^{\mathbf{r}})^{\sharp}_{s}(P),
       (p^{\mathbf{r}})^{\sharp}_{s_{0}}(Q_{0}),\ldots,
       (p^{\mathbf{r}})^{\sharp}_{s_{\bb{\mathbf{s}}-1}}(Q_{\bb{\mathbf{s}}-1})
      )
    \tag{by Lemma~\ref{L:aux}}
      \\
    &=
    \xi_{\mathbf{r},\mathbf{s},s}^{\mathbf{A}}
      ( \widehat{f}_{\mathbf{s},s}(P),
        \widehat{f}_{\mathbf{r},s_{0}}(Q_{0}),\ldots,
        \widehat{f}_{\mathbf{r},s_{\bb{\mathbf{s}}-1}}(Q_{\bb{\mathbf{s}}-1})
      ).
\end{align*}
\end{flushleft}
Therefore the $S^{\star}\times S$-sorted mapping $\widehat{f}$ is a homomorphism. Furthermore, $\widehat{f}\circ \eta^{\Sigma} = f$, because, for every $\mathbf{s}\in S^{\star}$, $s\in S$, and $\sigma\in \Sigma_{\mathbf{s},s}$, we have that%
\begin{align*}
\widehat{f}_{\mathbf{s},s}(\eta^{\Sigma}_{\mathbf{s},s}(\sigma)) &=
    (p^{\mathbf{s}})^{\sharp}_{s}(\sigma(v^{s_{0}}_{0},\ldots,v^{s_{\bb{\mathbf{s}}-1}}_{\bb{\mathbf{s}}-1}))
    \\&=
    \sigma^{\mathbf{A}_{\mathbf{s}}}
      (p^{\mathbf{s}}_{{s}_{0}}(v^{s_{0}}_{0}),
        \ldots,
       p^{\mathbf{s}}_{{s}_{\bb{\mathbf{s}}-1}}(v^{s_{\bb{\mathbf{s}}-1}}_{\bb{\mathbf{s}}-1})
      ) \\
    &=
    \xi_{\mathbf{s},\mathbf{s},s}^{\mathbf{A}}
      ( f_{(\mathbf{s},s)}(\sigma),
        (\pi^{\mathbf{s}}_{0})^{\mathbf{A}},
        \ldots,
        (\pi^{\mathbf{s}}_{\bb{\mathbf{s}}-1})^{\mathbf{A}}
      ) \\
    &=
    f_{\mathbf{s},s}(\sigma) \tag{by $\mathrm{H}_{2}$}.
\end{align*}
It is obvious that $\widehat{f}$ is the unique homomorphism such that $\widehat{f}\circ \eta^{\Sigma} = f$.  Henceforth $\mathbf{Tm}_{\mathrm{H}_{S}}(\Sigma)$ is isomorphic to $\mathbf{T}_{\mathrm{H}_{S}}(\Sigma)$.
\end{proof}

This isomorphism, together with the adjunction $\mathbf{T}_{\mathrm{H}_{S}}\dashv\mathrm{G}_{\mathrm{H}_{S}}$,  has as an immediate consequence that, for every $S$-sorted set $A$ and every $S$-sorted signature $\Sigma$, the sets $\mathrm{Hom}({\Sigma},\mathrm{G}_{\mathrm{H}_{S}}(\mathbf{O}_{\mathrm{H}_{S}}(A)))$, in 
$\mathsf{Set}^{S^{\star}\times S}$, and $\mathrm{Hom}(\mathbf{Tm}_{\mathrm{H}_{S}}(\Sigma),\mathbf{O}_{\mathrm{H}_{S}}(A))$, in $\mathsf{Alg}(\mathrm{H}_{S})$, are naturally isomorphic.

Actually, for an $S$-sorted set $A$ and an $S$-sorted signature $\Sigma$,
\begin{enumerate}
\item the mapping that assigns to a structure of $\Sigma$-algebra $F$ on $A$ (i.e., an $S^{\star}\times S$-sorted mapping $F$ from $\Sigma$ to $\mathrm{O}_{\mathrm{H}_{S}}(A)$) the unique  homomorphism of Hall algebras $\widehat{F} = (\widehat{F}_{w,s})_{(w,s)\in S^{\star}\times S}$ from $\mathbf{Tm}_{\mathrm{H}_{S}}(\Sigma)$ to $\mathbf{O}_{\mathrm{H}_{S}}(A)$ such that $\widehat{F}\circ \eta^{\Sigma} = F$, together with
\item the mapping that assigns to a homomorphism $h$ from $\mathbf{Tm}_{\mathrm{H}_{S}}(\Sigma)$ to $\mathbf{O}_{\mathrm{H}_{S}}(A)$, essentially, the algebraic structure $\mathrm{G}_{\mathrm{H}_{S}}(h)\circ \eta^{\Sigma}$ on $A$, where $\eta^{\Sigma}$ is the canonical embedding of $\Sigma$ into $\mathrm{T}_{\mathrm{H}_{S}}(\Sigma)$, 
\end{enumerate}
which are natural in both variables and mutually inverse bijections, provide the mentioned natural isomorphism.

\begin{remark}
For an $S$-sorted set $A$, a structure of $\Sigma$-algebra $F$ on $A$ and a word $\mathbf{s}\in S^{\star}$, the subfamily $(\widehat{F}_{\mathbf{s},s})_{s\in S}$ of $\widehat{F} = (\widehat{F}_{\mathbf{s},s})_{(\mathbf{s},s)\in S^{\star}\times S}$,
denoted by $\mathrm{d}^{(A,F)} = (\mathrm{d}^{(A,F)}_{\mathbf{s},s})_{s\in S}$, is the unique homomorphism from $\mathbf{T}_{\Sigma}(\vs{\mathbf{s}})$ to $(A,F)^{A_{\mathbf{s}}}$, the direct $A_{\mathbf{s}}$-power of $(A,F)$, such that $\mathrm{d}^{(A,F)}_{\mathbf{s}}\circ \eta^{\vs{\mathbf{s}}} = \mathrm{p}^{A}_{\mathbf{s}}$, where $\mathrm{p}^{A}_{\mathbf{s}}$ is the $S$-sorted mapping from $\vs{\mathbf{s}}$ to $A^{A_{\mathbf{s}}}$ defined, for every $s\in S$ and every $v^{s}_{i}\in(\vs{\mathbf{s}})_{s}$, as $\mathrm{p}^{A}_{\mathbf{s},s}(v^{s}_{i}) = \mathrm{pr}^{A}_{\mathbf{s},i}$.
\end{remark}

For a set of sorts $S$ and an $S$-sorted signature $\Sigma$, the derived operation symbols of $\Sigma$, i.e., the elements of the components $\mathrm{T}_{\Sigma}(\vs{\mathbf{s}})_{s}$ of $\mathrm{Tm}_{\mathrm{H}_{S}}(\Sigma)$ can be considered as the operation symbols of a new signature. The mappings $\varphi$ from a set of sorts $S$ to another $T$, together with the $S^{\star}\times S$-sorted mappings $d$ that assign to operation symbols of an $S$-sorted signature $\Sigma$ derived operations symbols of a $T$-sorted signature $\Lambda$, i.e., terms of the components $\mathrm{T}_{\Lambda}(\vs{\varphi^{\star}(\mathbf{s})})_{\varphi(s)}$ of $\mathrm{Tm}_{\mathrm{H}_{T}}(\Lambda)_{\varphi^{\star}\times\varphi}$, form a new class of morphisms between many-sorted signatures which we will call \emph{derivors}.

\begin{definition}\label{DDerivor}
Let $\mathbf{\Sigma} = (S,\Sigma)$ and $\mathbf{\Lambda} = (T,\Lambda)$ be many-sorted signatures. A \emph{derivor} from $\mathbf{\Sigma}$ to $\mathbf{\Lambda}$ is an ordered pair $\mathbf{d} = (\varphi,d)$, where $\varphi\colon S\mor T$ and  $d\colon\Sigma\mor\mathrm{Tm}_{\mathrm{H}_{T}}(\Lambda)_{\varphi^{\star}\times\varphi}$. Thus, for every $(\mathbf{s},s)\in S^{\star}\times S$, %
$$
d_{\mathbf{s},s}\colon \Sigma_{\mathbf{s},s}\mor\mathrm{Tm}_{\mathrm{H}_{T}}(\Lambda)_{\varphi^{\star}(\mathbf{s}),\varphi(s)}=
    \mathrm{T}_{\Lambda}(\vs{\varphi^{\star}(\mathbf{s})})_{\varphi(s)}.
$$
Hence, for every $\sigma\in \Sigma_{\mathbf{s},s}$, $d_{\mathbf{s},s}(\sigma)\in \mathrm{T}_{\Lambda}(\vs{\varphi^{\star}(\mathbf{s})})_{\varphi(s)}$, i.e., $d_{\mathbf{s},s}(\sigma)$ is a term of sort $\varphi(s)$ with variables in $\vs{\varphi^{\star}(\mathbf{s})}$.
\end{definition}

\begin{remark}\label{Rdop}
Let $\varphi$ be a mapping from the set of sorts $S$ to the set of sorts $T$ and let $V^{T}$ be a fixed $T$-sorted set with a countably infinity of variables in each component, say $V^{T}_{t} = \{v^{t}_{n}\mid n\in \mathbb{N}\}$ (which is isomorphic to $\{t\}\times\mathbb{N}$). Then, for every $\mathbf{s}\in S^{\star}$, we let $\vs{\varphi^{\star}(\mathbf{s})} = ((\vs{\varphi^{\star}(\mathbf{s})})_{t})_{t\in T}$ stand for the $T$-sorted set  defined, for every $t\in T$, as:
$$
(\vs{\varphi^{\star}(\mathbf{s})})_{t} = \{v^{t}_{0},\ldots,v^{t}_{\bb{\varphi^{\star}(\mathbf{s})}_{t}-1}\} (\cong \varphi^{\star}(\mathbf{s})^{-1}[\{t\}] = \{i\in \bb{\varphi^{\star}(\mathbf{s})}\mid \varphi(w_{i}) = t\}).
$$
Thus, for every $t\in T$, $(\vs{\varphi^{\star}(\mathbf{s})})_{t}$ is an initial segment of $V^{T}_{t}$ and 
$\vs{\varphi^{\star}(\mathbf{s})}$ is a finite subset of $V^{T}$. The latter is fulfilled because $\mathrm{supp}_{T}(\vs{\varphi^{\star}(\mathbf{s})}) = \mathrm{Im}(\varphi^{\star}(\mathbf{s}))$ and $\bb{\varphi^{\star}(\mathbf{s})} = \bb{\mathbf{s}}$ is finite. Moreover, $\mathrm{supp}_{T}(\vs{\varphi^{\star}(\mathbf{s})}) =\varphi[\mathrm{supp}_{S}(\vs{\mathbf{s}})]$. Recall  that $\varphi^{\star}(\mathbf{s}) = (\varphi(s_{i}))_{i\in \bb{\mathbf{s}}}$ and that, since $\varphi^{\star}(\lambda) = \lambda$, $\vs{\varphi^{\star}(\lambda)} = \varnothing^{T} = (\varnothing)_{t\in T}$.

If $\mathbf{d} = (\varphi,d)$ is a derivor from $\mathbf{\Sigma}$ to $\mathbf{\Lambda}$, $\mathbf{B}$ a 
$\mathbf{\Lambda}$-algebra, $(\mathbf{s},s)\in S^{\star}\times S$ and $\sigma\in \Sigma_{\mathbf{s},s}$, then the derived operation (or term operation or polynomial operation) of type $(\varphi^{\star}(\mathbf{s}),\varphi(s))$ determined by the term $d_{\mathbf{s},s}(\sigma)$ in $\mathbf{B}$, denoted by $d_{\mathbf{s},s}(\sigma)^{\mathbf{B}}$ or by $d(\sigma)^{\mathbf{B}}$ for short, is a mapping from $B_{\varphi^{\star}(\mathbf{s})}$ to $B_{\varphi(s)}$. In particular, for a $T$-sorted set $Y$, we have that the derived operation (or term operation or polynomial operation) of type $(\varphi^{\star}(\mathbf{s}),\varphi(s))$ determined by the term $d_{\mathbf{s},s}(\sigma)$ in $\mathbf{T}_{\Lambda}(Y)$, denoted by $d_{\mathbf{s},s}(\sigma)^{\mathbf{T}_{\Lambda}(Y)}$ or by  $d(\sigma)^{\mathbf{T}_{\Lambda}(Y)}$ for short, is a mapping from $\mathrm{T}_{\Lambda}(Y)_{\varphi^{\star}(\mathbf{s})}$ to $\mathrm{T}_{\Lambda}(Y)_{\varphi(s)}$. Notice that, for every $s\in S$ and every $\sigma\in \Sigma_{\lambda,s}$, $d_{\lambda,s}(\sigma)\in \mathrm{T}_{\Lambda}(\varnothing^{T})_{\varphi(s)}$, i.e., $d_{\lambda,s}(\sigma)$ is a closed term of sort $\varphi(s)$.  Therefore, for every $\mathbf{\Lambda}$-algebra $\mathbf{B}$, $d_{\lambda,s}(\sigma)^{\mathbf{B}}$ will be a mapping from $B_{\lambda}$ to $B_{\varphi(s)}$, i.e., essentially, a distinguished element of $B_{\varphi(s)}$. In particular, for a $T$-sorted set $Y$, $d_{\lambda,s}(\sigma)^{\mathbf{T}_{\Lambda}(Y)}$ will be a mapping from $\mathrm{T}_{\Lambda}(Y)_{\lambda}$ to $\mathrm{T}_{\Lambda}(Y)_{\varphi(s)}$, i.e., essentially, a distinguished element of $\mathrm{T}_{\Lambda}(Y)_{\varphi(s)}$. Let us note that $d(\sigma)^{\mathbf{B}}$ is the value of $\mathrm{d}^{\mathbf{B}}_{\varphi^{\star}(\mathbf{s}),s}$ at $d_{\mathbf{s},s}(\sigma)$ and that $d(\sigma)^{\mathbf{T}_{\Lambda}(Y)}$ is the value of $\mathrm{d}^{\mathbf{T}_{\Lambda}(Y)}_{\varphi^{\star}(\mathbf{s}),s}$ at $d_{\mathbf{s},s}(\sigma)$, where 
$\mathrm{d}^{\mathbf{B}}_{\varphi^{\star}(\mathbf{s})}$ is the unique homomorphism from $\mathbf{T}_{\Lambda}(\vs{\varphi^{\star}(\mathbf{s})})$ to $\mathbf{B}^{B_{\varphi^{\star}(\mathbf{s})}}$ for which the left diagram below commutes, and $\mathrm{d}^{\mathbf{T}_{\Lambda}(Y)}_{\varphi^{\star}(\mathbf{s})}$ the unique homomorphism from $\mathbf{T}_{\Lambda}(\vs{\varphi^{\star}(\mathbf{s})})$ to $\mathbf{T}_{\Lambda}(Y)^{\mathrm{T}_{\Lambda}(Y)_{\varphi^{\star}(\mathbf{s})}}$ for which the right diagram below commutes.
$$
\begin{aligned}
\xymatrix@C=40pt@R=33pt{
\vs{\varphi^{\star}(\mathbf{s})}
  \ar[r]^-{\eta^{\vs{\varphi^{\star}(\mathbf{s})}}}
  \ar[rd]_-{\mathrm{p}^{B}_{\varphi^{\star}(\mathbf{s})}} &
\mathrm{T}_{\Lambda}(\vs{\varphi^{\star}(\mathbf{s})}) \ar[d]^{\mathrm{d}^{\mathbf{B}}_{\varphi^{\star}(\mathbf{s})}}
  \\
& B^{B_{\varphi^{\star}(\mathbf{s})}}
}
\end{aligned}
\quad \text{} \qquad\;
\begin{aligned}
\xymatrix@C=40pt@R=33pt{
\vs{\varphi^{\star}(\mathbf{s})}
  \ar[r]^-{\eta^{\vs{\varphi^{\star}(\mathbf{s})}}}
  \ar[rd]_-{\mathrm{p}^{\mathrm{T}_{\Lambda}(Y)}_{\varphi^{\star}(\mathbf{s})}} &
\mathrm{T}_{\Lambda}(\vs{\varphi^{\star}(\mathbf{s})}) \ar[d]^{\mathrm{d}^{\mathbf{T}_{\Lambda}(Y)}_{\varphi^{\star}(\mathbf{s})}}
  \\
& \mathrm{T}_{\Lambda}(Y)^{\mathrm{T}_{\Lambda}(Y)_{\varphi^{\star}(\mathbf{s})}}
}
\end{aligned}
$$
Recall that, for every $t\in T$ and every $i\in (\vs{\varphi^{\star}(\mathbf{s})})_{t}$ (thus $(\varphi^{\star}(\mathbf{s}))(i) = \varphi(s_{i}) = t$), $\mathrm{p}^{B}_{\varphi^{\star}(\mathbf{s}),t}(i)$ is $\mathrm{pr}^{B}_{\varphi^{\star}(\mathbf{s}),i}$, i.e., the $\varphi^{\star}(\mathbf{s})$-ary, $i$-th canonical projection from $B_{\varphi^{\star}(\mathbf{s})}$ to $B_{t} = B_{(\varphi^{\star}(\mathbf{s}))(i)}$.
\end{remark}

\begin{remark}
While derivors are not the most general type of morphism that might be considered between many-sorted signatures---for instance, one could consider polyderivors, see~\cite{CVST10b}---, they are an important class of such morphisms. One reason for its relevance is its formal properties (see below), another that there are many mathematical examples of them which are of interest (see~\cite{CVST10b}).
\end{remark}

Since by, Proposition~\ref{iso:FrH-TerH}, $\mathbf{Tm}_{\mathrm{H}_{T}}(\Lambda)$ is isomorphic to 
$\mathbf{T}_{\mathrm{H}_{T}}(\Lambda)$, the derivors can be defined, alternative, but equivalently, as ordered pairs  $\mathbf{d} = (\varphi,d)$ with $d\colon\Sigma\mor \mathrm{T}_{\mathrm{H}_{T}}(\Lambda)$.

On the other hand, every mapping $\varphi\colon S\mor T$ determines a functor from the category  $\mathsf{Alg}(\mathrm{H}_{T})$ to the category $\mathsf{Alg}(\mathrm{H}_{S})$, so $\mathrm{Tm}_{\mathrm{H}_{T}}(\Lambda)_{\varphi^{\star}\times\varphi}$ is in its turn e\-quipped with a structure of Hall algebra for $S$, which will allow us, in particular, to define the composition of derivors. We next show the existence of such a functor by defining a morphism of algebraic presentations from  $(\Sigma^{\mathrm{H}_{S}},\mathcal{E}^{\mathrm{H}_{S}})$ to $(\Sigma^{\mathrm{H}_{T}},\mathcal{E}^{\mathrm{H}_{T}})$.

\begin{proposition}
Let $\varphi\colon S\mor T$ be a mapping between sets of sorts and $h^{\varphi}$ the $S^{\star}\times S$-sorted mapping from $\Sigma^{\mathrm{H}_{S}}$ to $\Sigma^{\mathrm{H}_{T}}_{\varphi^{\star}\times \varphi}$ defined as follows: 
\begin{enumerate}
\item For every $\mathbf{s}\in S^{\star}$ and every $i\in\bb{\mathbf{s}}$, $h^{\varphi}(\pi^{\mathbf{s}}_{i}) = \pi^{\varphi^{\star}(\mathbf{s})}_{i}$.

\item  For every $\mathbf{r},\,\mathbf{s}\in S^{\star}$ and every $s\in S$,
  $h^{\varphi}(\xi_{\mathbf{r},\mathbf{s},s}) = \xi_{\varphi^{\star}(\mathbf{r}),\varphi^{\star}(\mathbf{s}),\varphi(s)}$.
\end{enumerate}
Then $(\varphi^{\star}\times \varphi, h^{\varphi}) \colon(S^{\star}\times S,\Sigma^{\mathrm{H}_{S}},\mathcal{E}^{\mathrm{H}_{S}})\mor
(T^{\star}\times T,\Sigma^{\mathrm{H}_{T}},\mathcal{E}^{\mathrm{H}_{T}})$ is a morphism of algebraic presentations.
\end{proposition}

The morphisms of algebraic presentations determine functors in the opposite direction between the associated categories of algebras. Therefore, each mapping $\varphi\colon S\mor T$ between sets of sorts, determines a functor $(\varphi^{\star}\times \varphi, h^{\varphi})^{\ast}$ from $\mathsf{Alg}(\mathrm{H}_{T})$ to $\mathsf{Alg}(\mathrm{H}_{S})$, which transforms Hall algebras for $T$ into Hall algebras for $S$. The action of the functor on the free Hall algebra on a $T$-sorted signature $\Lambda$ is a Hall algebra for $S$, whose underlying $S^{\star}\times S$-sorted set is
$\mathrm{Tm}_{\mathrm{H}_{T}}(\Lambda)_{\varphi^{\star}\times \varphi}$.

If $\mathbf{d}\colon\mathbf{\Sigma}\mor\mathbf{\Lambda}$ is a derivor, then
$d\colon \Sigma\mor \mathrm{Tm}_{\mathrm{H}_{T}}(\Lambda)_{\varphi^{\star}\times\varphi}$ determines a homomorphism of Hall algebras $d^{\sharp}\colon\mathbf{Tm}_{\mathrm{H}_{S}}(\Sigma)\mor
\mathbf{Tm}_{\mathrm{H}_{T}}(\Lambda)_{\varphi^{\star}\times\varphi}$. Thus, for every $(\mathbf{s},s)\in S^{\star}\times S$, $d^{\sharp}_{\mathbf{s},s}$ sends terms in $\mathrm{T}_{\Sigma}(\vs{\mathbf{s}})_{s}$ to terms in $\mathrm{T}_{\Lambda}(\vs{\varphi^{\star}(\mathbf{s})})_{\varphi(s)}$.

\begin{definition}\label{DDerCompId}
Let $\mathbf{d}\colon\mathbf{\Sigma}\mor\mathbf{\Lambda}$ and $\mathbf{e}\colon\mathbf{\Lambda}\mor\mathbf{\Omega}$ be  derivors. Then $\mathbf{e}\circ \mathbf{d} = (\psi,e)\circ(\varphi,d)$, the \emph{composition} of $\mathbf{d}$ and  $\mathbf{e}$, is the derivor $(\psi\circ\varphi,e^{\sharp}_{\varphi^{\star}\times
\varphi}\circ d)$, where $e^{\sharp}_{\varphi^{\star}\times\varphi}\circ d$ is obtained from
$$
\begin{aligned}
\xymatrix@C=40pt@R=30pt{
\Lambda
  \ar[r]^-{\eta^{\Lambda}}
  \ar[rd]_-{e} &
\mathrm{Tm}_{\mathrm{H}_{T}}(\Lambda) \ar[d]^{e^{\sharp}}
  \\
& \mathrm{Tm}_{\mathrm{H}_{U}}(\Omega)_{\psi^{\star}\times \psi}
}
\end{aligned}
\qquad \text{as} \qquad\;\;
\begin{aligned}
\xymatrix@C=40pt@R=30pt{
\mathrm{Tm}_{\mathrm{H}_{T}}(\Lambda)_{\varphi^{\star}\times \varphi}
\ar[d]^{e^{\sharp}_{\varphi^{\star}\times \varphi}} &
\Sigma \ar[l]_-{d} \\
{\mathrm{Tm}_{\mathrm{H}_{U}}(\Omega)_{\psi^{\star}\times\psi}}_{\varphi^{\star}\times\varphi}
}
\end{aligned}
$$
being $e^{\sharp}$ the canonical extension of $e$ to the free Hall algebra on $\Lambda$.

For every many-sorted signature $\mathbf{\Sigma} = (S,\Sigma)$, the \emph{identity} at $(S,\Sigma)$
is $(\id^{S},\eta^{\Sigma})$.
\end{definition}

The just stated definition allows us to form a category of many-sorted signatures whose morphisms are the derivors.

\begin{proposition}\label{PSigdCat}
The many-sorted signatures together with the derivors constitute a category, denoted by $\mathsf{Sig}_{\mathfrak{d}}$.
\end{proposition}

\begin{remark}
Let $\mathsf{Sig}$ be the category whose objects are the many-sorted signatures and whose morphisms from $\mathbf{\Sigma}$ to $\mathbf{\Lambda}$ are the ordered pairs $(\varphi,d)$, where $\varphi$ is a mapping from $S$ to $T$ and $d$ an $S^{\star}\times S$-sorted mapping from $\Sigma$ to $\Lambda_{\varphi^{\star}\times \varphi}$ (thus a mapping in $\mathsf{Sig}(S)$). Note that $\mathsf{Sig}$ is the category obtained by means of the Grothendieck construction for an  appropriate contravariant functor from $\mathsf{Set}$ to $\mathsf{Cat}$.
We next show that the category $\mathsf{Sig}_{\mathfrak{d}}$ can be obtained as the Kleisli category for a suitable monad. In fact, for every  set of sorts $S$, we have the adjoint situation $\mathbf{T}_{\mathrm{H}_{S}}\dashv\mathrm{G}_{\mathrm{H}_{S}}$
and, thus, a monad on $\mathsf{Sig}(S)$ which we will denote by $\mathbb{T}_{\mathrm{H}_{S}} = (\mathrm{T}_{\mathrm{H}_{S}},\eta^{\mathrm{H}_{S}},\mu^{\mathrm{H}_{S}})$. Then the ordered triple $\boldsymbol{\mathfrak{d}} = (\mathfrak{d},\eta,\mu)$ in which (1) $\mathfrak{d}$ is the endofunctor at $\mathsf{Sig}$ that sends $(S,\Sigma)$ to $(S,\mathrm{T}_{\mathrm{H}_{S}}(\Sigma))$ and a morphism $(\varphi,d)$ from $(S,\Sigma)$ to  $(T,\Lambda)$ to the morphism $(\varphi,d^{\sharp})$ from $(S,\mathrm{T}_{\mathrm{H}_{S}}(\Sigma))$ to  $(T,\mathrm{T}_{\mathrm{H}_{T}}(\Lambda))$, (2) $\eta$ the natural transformation from $\mathrm{Id}_{\mathsf{Sig}}$ to  $\mathfrak{d}$ that sends $(S,\Sigma)$ to $\eta_{(S,\Sigma)} = (\mathrm{id},\eta^{\mathrm{H}_{S}}_{\Sigma})$, and (3) $\mu$ the natural transformation from $\mathfrak{d}\circ \mathfrak{d}$ to $\mathfrak{d}$ that sends $(S,\Sigma)$ to $\mu_{(S,\Sigma)} = (\mathrm{id},\mu^{\mathrm{H}_{S}}_{\Sigma})$, is a monad on $\mathsf{Sig}$ and $\mathsf{Kl}(\boldsymbol{\mathfrak{d}})$, the Kleisli category for $\boldsymbol{\mathfrak{d}}$, is isomorphic to $\mathsf{Sig}_{\mathfrak{d}}$. This shows, in our opinion, the mathematical naturalness of the notion of derivor. Moreover, by defining a suitable notion of transformation between derivors one can equip the category $\mathsf{Sig}_{\mathfrak{d}}$ with a structure of $2$-category (this was, in fact, already done in~\cite{CVST10b} for polyderivors).

\end{remark}

\begin{remark}
Since $\mathsf{Sig}$ has coproducts, $\mathsf{Kl}(\boldsymbol{\mathfrak{d}})\cong \mathsf{Sig}_{\mathfrak{d}}$ has coproducts. In fact, let $(\mathbf{\Sigma}^{i})_{i\in I}$ be an $I$-indexed family in $\mathrm{Ob}(\mathsf{Kl}(\boldsymbol{\mathfrak{d}}))$ and $(\coprod_{i\in I}\mathbf{\Sigma}^{i},(\mathrm{in}^{\mathbf{\Sigma}^{i}})_{i\in I})$ the coproduct of $(\mathbf{\Sigma}^{i})_{i\in I}$ in $\mathsf{Sig}$, where  the first coordinate of $\coprod_{i\in I}\mathbf{\Sigma}^{i}$ is $\coprod_{i\in I}S_{i}$ and the second coordinate of $\coprod_{i\in I}\mathbf{\Sigma}^{i}$, denoted by $\coprod_{i\in I}\mathrm{\Sigma}^{i}$, is defined, for every $(((s_{i_{j}},i_{j})_{j\in n}),(s,i))\in (\coprod_{i\in I}S_{i})^{\star}\times \coprod_{i\in I}S_{i}$, as follows: 
\begin{equation*}
\textstyle
(\coprod_{i\in I}\mathrm{\Sigma}^{i})_{(((s_{i_{j}},i_{j})_{j\in n}),(s,i))} = 
\begin{cases}
\Sigma^{i}_{\mathbf{s},s}, \text{\,\,if $\exists\,\mathbf{s}\in S_{i}^{\star}$ such that $(s_{i_{j}},i_{j})_{j\in n} = \mathrm{in}_{i}^{\star}(\mathbf{s})$;}\\
\varnothing, \text{\,\,otherwise}.
\end{cases}
\end{equation*}
Let us recall that $\mathrm{in}_{i}^{\star}$ is the monoid homomorphism from the free monoid on $S_{i}$ to the free monoid on $\coprod_{i\in I}S_{i}$ that sends $\mathbf{s} = (s_{j})_{j\in n}\in S_{i}^{\star}$ to 
$(s_{j},i)_{j\in n}\in (\coprod_{i\in I}S_{i})^{\star}$. Therefore, 
$(s_{i_{j}},i_{j})_{j\in n} = \mathrm{in}_{i}^{\star}(\mathbf{s})$ if and only if, for every $j\in n$, $i_{j} = i$ and $s_{i_{j}} = w_{j}$.
Then $(\boldsymbol{\mathfrak{d}}(\coprod_{i\in I}\mathbf{\Sigma}^{i}),(\eta_{\coprod_{i\in I}\mathbf{\Sigma}^{i}}\circ \mathrm{in}^{\mathbf{\Sigma}^{i}})_{i\in I})$ will be the required coproduct of $(\mathbf{\Sigma}^{i})_{i\in I}$ in $\mathsf{Kl}(\boldsymbol{\mathfrak{d}})$. This property of the category $\mathsf{Sig}_{\mathfrak{d}}$ will be essential in our work because to be able to compare Curry-Howard towers it will be necessary, as an intermediate step, to construct, by means of the formation of coproducts, from given many-sorted  signatures $(\mathbf{\Sigma}^{i})_{i\in I}$, a new many-sorted signature $\coprod_{i\in I}\mathbf{\Sigma}^{i}$ (finally, from $\coprod_{i\in I}\mathbf{\Sigma}^{i}$, we will construct the many-sorted signature $\biguplus_{i\in I}\mathbf{\Sigma}^{i}$ which will be the one of interest for us).

\end{remark}

We next associate to every derivor $\mathbf{d}\colon\mathbf{\Sigma}\mor\mathbf{\Lambda}$ a functor from  $\mathsf{Alg}(\mathbf{\Lambda}) (= \mathsf{Alg}(\Lambda))$ to $\mathsf{Alg}(\mathbf{\Sigma}) (= \mathsf{Alg}(\Sigma))$.

\begin{proposition}\label{PFunSig}
Let $\mathbf{d}\colon\mathbf{\Sigma}\mor\mathbf{\Lambda}$ be a morphism in
$\mathsf{Sig}_{\mathfrak{d}}$. Then $\mathrm{Alg}_{\mathfrak{d}}(\mathbf{d})$, also denoted by $\mathbf{d}^{\ast}_{\mathfrak{d}}$ for short, is the functor from  $\mathsf{Alg}(\mathbf{\Lambda})$ to $\mathsf{Alg}(\mathbf{\Sigma})$ that sends $(B,G)$ to $(B_{\varphi},G^{\mathbf{d}})$ and a homomorphism $f$ from $(B,G)$ to $(B',G')$ to the homomorphism $f_{\varphi}$ from $(B_{\varphi},G^{\mathbf{d}})$ to  $(B'_{\varphi},{G'}^{\mathbf{d}})$,
where, for every $\mathbf{\Lambda}$-algebra $(B,G)$,
$G^{\mathbf{d}} = G^{\sharp}_{\varphi^{\star}\times \varphi}\circ d$
is obtained from
$$
\begin{aligned}
\xymatrix@C=40pt@R=33pt{
\Lambda
  \ar[r]^-{\eta_{\Lambda}}
  \ar[rd]_-{G} &
\mathrm{Tm}_{\mathrm{H}_{T}}(\Lambda) \ar[d]^{G^{\sharp}}
  \\
& \mathrm{O}_{\mathrm{H}_{T}}(B)
}
\end{aligned}
\qquad\quad \text{as} \quad
\begin{aligned}
\xymatrix@C=40pt@R=33pt{
\mathrm{Tm}_{\mathrm{H}_{T}}(\Lambda)_{\varphi^{\star}\times \varphi}
\ar[d]^{G^{\sharp}_{\varphi^{\star}\times \varphi}} &
\Sigma \ar[l]_-{d}\\
{\mathrm{O}_{\mathrm{H}_{T}}(B)_{\varphi^{\star}\times\varphi}}=
\mathrm{O}_{\mathrm{H}_{S}}(B_{\varphi})
}
\end{aligned}
$$
\end{proposition}

\begin{proof}
For every $\mathbf{\Lambda}$-algebra $(B,G)$, $G\colon\Lambda\mor\mathrm{O}_{\mathrm{H}_{T}}(B)$, and since  $\mathbf{O}_{\mathrm{H}_{T}}(B)$ is a Hall algebra, $G$  can be extended up to the free Hall algebra on $\Lambda$. Moreover, we have that $\mathrm{O}_{\mathrm{H}_{T}}(B)_{\varphi^{\star}\times \varphi} = \mathrm{O}_{\mathrm{H}_{S}}(B_{\varphi})$ since, for every $(\mathbf{s},s)\in S^{\star}\times S$ it holds that
\allowdisplaybreaks
\begin{align*}
    (\mathrm{O}_{\mathrm{H}_{T}}(B)_{\varphi^{\star}\times\varphi})_{\mathbf{s},s} &=
    \mathrm{O}_{\mathrm{H}_{T}}(B)_{\varphi^{\star}(\mathbf{s}),\varphi(s)} \\
    &= B_{\varphi^{\star}(\mathbf{s})}\mor B_{\varphi(s)} \\
    &= B_{(\varphi(s_{0}),\ldots,\varphi(s_{\bb{\mathbf{s}}-1}))}\mor B_{\varphi(s)} \\
    &= \textstyle\prod(B_{\varphi(s_{i})}\mid i\in\bb{\mathbf{s}}) \mor B_{\varphi(s)} \\
    &= \textstyle\prod((B_{\varphi})_{s_{i}}\mid i\in\bb{\mathbf{s}}) \mor (B_{\varphi})_{s} \\
    &= (B_{\varphi})_{\mathbf{s}} \mor (B_{\varphi})_{s} \\
    &= \mathrm{O}_{\mathrm{H}_{S}}(B_{\varphi})_{\mathbf{s},s},
\end{align*}
thus, so defined, $G^{(\varphi,d)}$ is  an algebraic structure on $B_{\varphi}$.

Let $f\colon(B,G)\mor(B',G')$ be a homomorphism of $\mathbf{\Lambda}$-algebras, $(\mathbf{s},s)\in S^{\star}\times S$, and $\sigma\in\Sigma_{\mathbf{s},s}$. Then  $f_{\varphi}\colon(B_{\varphi},G^{\mathbf{d}}) \mor (B'_{\varphi},{G'}^{\mathbf{d}})$ is a homomorphism of $\mathbf{\Sigma}$-algebras, because $G^{\mathbf{d}}(\sigma)$ is a term operation and, consequently, the following diagram
\begin{center}
\begin{tikzpicture}
[ACliment/.style={-{To [angle'=45, length=5.75pt, width=4pt, round]}},
RHACliment/.style={right hook-{To [angle'=45, length=5.75pt, width=4pt, round]}, font=\scriptsize}, 
scale=1
]

\node[] (s) at (0,0) {${B_{\varphi}}_{\mathbf{s}}$};
\node[] (s') at (0,-2) {${B'_{\varphi}}_{\mathbf{s}}$};
\node[] (Vs) at (4,0) {$V^{S}_{\mathbf{s}}$};
\node[] (Vs') at (4,-2) {$V^{S}_{\mathbf{s}'}$};

\draw[ACliment]  (s) to node [above ] {$G^{(\varphi,d)}(\sigma)$} (Vs);
\draw[ACliment]  (s') to node [below ] {$G'^{(\varphi,d)}(\sigma)$} (Vs');
\draw[ACliment]  (s) to node [left ] {${f_{\varphi}}_{\mathbf{s}}$} (s');
\draw[ACliment]  (Vs) to node [right ] {${f_{\varphi}}_{s}$} (Vs');
\end{tikzpicture}
\end{center}
commutes. Moreover, $(g\circ f)_{\varphi}=g_{\varphi}\circ f_{\varphi}$, so that  $\mathrm{Alg}_{\mathfrak{d}}(\mathbf{d})$ is a functor.
\end{proof}

From the definition of the functor $\mathrm{Alg}_{\mathfrak{d}}$, for every derivor $\mathbf{d}\colon\mathbf{\Sigma}\mor\mathbf{\Lambda}$, it is obvious that the following diagram %
\begin{center}
\begin{tikzpicture}
[ACliment/.style={-{To [angle'=45, length=5.75pt, width=4pt, round]}},
RHACliment/.style={right hook-{To [angle'=45, length=5.75pt, width=4pt, round]}, font=\scriptsize}, 
scale=1
]

\node[] (s) at (0,0) {$\mathsf{Alg}(\mathbf{\Sigma}$};
\node[] (s') at (0,-2) {$\mathsf{Alg}(\mathbf{\Lambda})$};
\node[] (Vs) at (4,0) {$\mathsf{Set}^{S}$};
\node[] (Vs') at (4,-2) {$\mathsf{Set}^{T}$};

\draw[ACliment]  (s) to node [above ] {$\G_{\mathbf{\Sigma}}$} (Vs);
\draw[ACliment]  (s') to node [below ] {$\G_{\mathbf{\Lambda}}$} (Vs');
\draw[ACliment]  (s') to node [left ] {$\mathrm{Alg}_{\mathfrak{d}}(\mathbf{d})$} (s);
\draw[ACliment]  (Vs') to node [right ] {$\Delta_{\varphi}$} (Vs);
\end{tikzpicture}
\end{center}
commutes.

The previous construction can be extended up to a contravariant functor from the category $\mathsf{Sig}_{\mathfrak{d}}$ to the category $\mathsf{Cat}$.

\begin{proposition}\label{PAlgFun}
From $\mathsf{Sig}_{\mathfrak{d}}$ to $\mathsf{Cat}$ there exists a contravariant functor,
denoted by $\mathrm{Alg}_{\mathfrak{d}}$, that sends $(S,\Sigma)$ to $\mathsf{Alg}(S,\Sigma)$ and a morphism  $(\varphi,d)$ from $(S,\Sigma)$ to $(T,\Lambda)$ to the functor $\mathrm{Alg}_{\mathfrak{d}}(\varphi,d)$ from  $\mathsf{Alg}(T,\Lambda)$ to $\mathsf{Alg}(S,\Sigma)$.
\end{proposition}

\begin{proof}
Given $(\varphi,d)\colon(S,\Sigma) \mor (T,\Lambda)$ and $(\psi,e)\colon (T,\Lambda)\mor (U,\Omega)$, we show that %
$\mathrm{Alg}_{\mathfrak{d}}(\varphi,d)\circ\mathrm{Alg}_{\mathfrak{d}}(\psi,e)= \mathrm{Alg}_{\mathfrak{d}}((\psi,e)\circ(\varphi,d))$.

Let $(A,F)$ be a $(U,\Omega)$-algebra. Then ${A_{\psi}}_{\varphi}=A_{\psi\circ\varphi}$. Moreover, we have that %
\begin{align*}
{F^{(\psi,e)}}^{(\varphi,d)}
  &=
(F^{\sharp}_{\psi^{\star}\times\psi}\circ e)^{(\varphi,d)} \\
  &=
(F^{\sharp}_{\psi^{\star}\times\psi}\circ e)^{\sharp}_{\varphi^{\star}\times\varphi}\circ d \\
  &=
({F^{\sharp}_{\psi^{\star}\times\psi}}_{\varphi^{\star}\times\varphi}\circ
    e^{\sharp}_{\varphi^{\star}\times\varphi})\circ d \\
  &=
{F^{\sharp}_{\psi^{\star}\times\psi}}_{\varphi^{\star}\times\varphi}\circ
    (e^{\sharp}_{\varphi^{\star}\times\varphi}\circ d )\\
  &=
F^{\sharp}_{(\psi\circ\varphi)^{\sharp}\times(\psi\circ\varphi)}
    \circ (e^{\sharp}_{\varphi^{\star}\times\varphi}\circ d )\\
  &=
F^{((\psi\circ\varphi),e^{\sharp}_{\varphi^{\star}\times\varphi}\circ d)} \\
  &=
F^{(\psi,e)\circ(\varphi,d)}
\end{align*}
thus $({A_{\psi}}_{\varphi},{F^{(\psi,e)}}^{(\varphi,d)})= (A_{\psi\circ\varphi},F^{(\psi,e)\circ (\varphi,d)})$.  Finally, if $f$ is a homomorphism of $(U,\Omega)$-algebras, then ${f_{\psi}}_{\varphi} = f_{\psi\circ\varphi}$.
\end{proof}

\begin{definition}
The category $\mathsf{Alg}_{\mathfrak{d}}$ is $\int^{\mathsf{Sig}_{\scriptscriptstyle\mathfrak{d}}}\mathrm{Alg}_{\mathfrak{d}}$, i.e., the category obtained by means of the Grothendieck construction applied to the contravariant functor $\mathrm{Alg}_{\mathfrak{d}}$. Thus the category $\mathsf{Alg}_{\mathfrak{d}}$ has as objects the ordered pairs $(\mathbf{\Sigma},\mathbf{A})$, with $\mathbf{\Sigma}$ a many-sorted signature and $\mathbf{A}$ a $\Sigma$-algebra, and as morphisms from $(\mathbf{\Sigma},\mathbf{A})$ to $(\mathbf{\Lambda},\mathbf{B})$ the ordered pairs $(\mathbf{d},f)$, with $\mathbf{d}$ a derivor from $\mathbf{\Sigma}$ to $\mathbf{\Lambda}$ and $f$ a homomorphism of $\Sigma$-algebras from $(A,F)$ to $(B_{\varphi},G^{\mathbf{d}})$.

Moreover, we let $\mathsf{FAlg}_{\mathfrak{d}}$ stand for the full subcategory of $\mathsf{Alg}_{\mathfrak{d}}$ whose objects are the ordered pairs $(\mathbf{\Sigma},\mathbf{T}_{\Sigma}(X))$, with 
$\mathbf{\Sigma} = (S,\Sigma)$ a many-sorted signature and $X$ an $S$-sorted set.
\end{definition}

\begin{remark}
There exists a pseudo-functor $\mathrm{Alg}_{\boldsymbol{\mathfrak{d}}}$ from the 2-category $\mathsf{Sig}_{\boldsymbol{\mathfrak{d}}}$ to the 2-category $\mathsf{Cat}$, contravariant in the 1-cells, i.e., the derivors, and covariant in the 2-cells, i.e., the transformations between derivors (this was already done in~\cite{CVST10b} for polyderivors).
\end{remark}

\section{zeroth-order rewriting systems}

In this section we define the concept of zeroth-order many-sorted rewriting system and, based on the notion of derivor, that of morphism between zeroth-order many-sorted rewriting systems. We show that they constitute a category denoted by $\mathsf{Rws}_{\mathfrak{d}}^{(0)}$. Furthermore, we define the zeroth-order tower associated to a zeroth-order many-sorted rewriting system and morphisms between zeroth-order towers. We show that they constitute a category denoted by $\mathsf{Tw}_{\mathfrak{d}}^{(0)}$. Finally, we show that both categories are isomorphic.

\begin{definition}\label{DRws0Mor}
A \emph{zeroth-order many-sorted rewriting system} or simply a \emph{zeroth-order rewriting system} is a triple $(S, \Sigma, X)$, often denoted by $\boldsymbol{\mathcal{A}}^{(0)}$ for short, where $(S, \Sigma)$ is a many-sorted signature, i.e., $S$ is a set of sorts and $\Sigma$ an $S$-sorted signature, and $X$ an $S$-sorted set. We will also call $(S,\Sigma,X)$ a \emph{zeroth-order $S$-sorted rewriting system}.

A \emph{morphism} of zeroth-order many-sorted rewriting systems, or simply a \emph{zeroth-order morphism}, from $\boldsymbol{\mathcal{A}}^{(0)} = (S, \Sigma, X)$ to $\boldsymbol{\mathcal{B}}^{(0)} = (T, \Lambda, Y)$ is an ordered triple $(\boldsymbol{\mathcal{A}}^{(0)}, \mathbf{f}^{(0)}, \boldsymbol{\mathcal{B}}^{(0)})$, denoted by $\mathbf{f}^{(0)} \colon \boldsymbol{\mathcal{A}}^{(0)} \mor \boldsymbol{\mathcal{B}}^{(0)}$ for short, in which $\mathbf{f}^{(0)} = (\varphi, c, f^{(0)})$ is the ordered triple where
\begin{enumerate}
\item
$\mathbf{c} = (\varphi, c)$, the \emph{underlying derivor} of $\mathbf{f}^{(0)}$, is a derivor from $(S, \Sigma)$ to $(T,\Lambda)$, as introduced in Definition~\ref{DDerivor}, and
\item
$f^{(0)} \colon X \mor \mathrm{T}_{\Lambda}(Y)_{\varphi}$ is an $S$-sorted mapping where 
$\T_{\Lambda}(Y)_{\varphi}$ is the underlying $S$-sorted set of the $\Sigma$-algebra $\mathbf{c}^{\ast}_{\mathfrak{d}}(\mathbf{T}_{\Lambda}(Y))$.
\end{enumerate}
The alternative notation $\mathbf{f}^{(0)} = (\mathbf{c}, f^{(0)})$ will also be used.
\end{definition}

We next define the identity morphism at a zeroth-order many-sorted rewriting system and the composite morphism between zeroth-order many-sorted rewriting systems. 

\begin{definition}
\label{DIdRws0}
Let $\boldsymbol{\mathcal{A}}^{(0)} = (S, \Sigma, X)$ be a zeroth-order many-sorted rewriting system. The identity morphism at $\boldsymbol{\mathcal{A}}^{(0)}$, denoted by $\mathrm{id}^{\boldsymbol{\mathcal{A}}^{(0)}}$, is given by
$$
\mathrm{id}^{\boldsymbol{\mathcal{A}}^{(0)}} = (\mathrm{id}^{\mathbf{\Sigma}}, \eta^{X}).
$$
where $\mathrm{id}^{\mathbf{\Sigma}}=(\mathrm{id}^{S}, \eta^{\Sigma})$ is the identity derivor at the many-sorted signature $\mathbf{\Sigma}=(S,\Sigma)$ introduced at Definition~\ref{DDerCompId} and $\eta^{X}$ is the standard inclusion of generators from $X$ to $\T_{\Sigma}(X)$.
\end{definition}

\begin{definition}
\label{DCompRws0}
Let $\mathbf{f}^{(0)} = (\mathbf{c}, f^{(0)})$ be a zeroth-order morphism from $\boldsymbol{\mathcal{A}}^{(0)}$ to $\boldsymbol{\mathcal{B}}^{(0)}$ and $\mathbf{g}^{(0)} = (\mathbf{d}, g^{(0)})$ a zeroth-order morphism from $\boldsymbol{\mathcal{B}}^{(0)}$ to $\boldsymbol{\mathcal{C}}^{(0)}$. Its zeroth-order \emph{composition} morphism, from $\boldsymbol{\mathcal{A}}^{(0)}$ to $\boldsymbol{\mathcal{C}}^{(0)}$, is given by
$$
\mathbf{g}^{(0)} \circ \mathbf{f}^{(0)}
=
\left(\mathbf{d}, g^{(0)}\right) \circ \left(\mathbf{c}, f^{(0)}\right)
=
\left(
\mathbf{d} \circ \mathbf{c},
g^{(0)\sharp}_{\varphi} \circ f^{(0)}
\right).
$$
Note that $\mathbf{d} \circ \mathbf{c}$ is the composition derivor, introduced at Definition~\ref{DDerCompId} and the construction of the mapping $g^{(0)\sharp}_{\varphi} \circ f^{(0)}$ is depicted in Figure~\ref{FCompRws0} being $g^{(0)\sharp}$ the canonical $\Lambda$-homomorphism extension of $g^{(0)}$ to the free $\Lambda$-algebra on $Y$.

\begin{figure}[t]
$$
\xymatrix{
Y
  \ar[r]^-{\eta^Y}
  \ar[rd]_-{g^{(0)}}
&
\mathrm{T}_{\Lambda}(Y)
  \ar[d]^-{g^{(0)\sharp}}
&&
\mathrm{T}_\Lambda(Y)_{\varphi}
  \ar[d]_-{g^{(0)\sharp}_{\varphi}} 
&
X
  \ar[l]_-{f^{(0)}}
\\
& 
\mathrm{T}_{\Omega}(Z)_{\psi}
&&
(\mathrm{T}_{\Omega}(Z)_{\psi})_{\varphi}
\\
}
$$
\caption{Construction of the mapping $g^{(0)\sharp}_{\varphi} \circ f^{(0)}$}
\label{FCompRws0}
\end{figure}
\end{definition}

\begin{remark}
\label{RCompExt}
Let $\mathbf{f}^{(0)} = (\mathbf{c}, f^{(0)})$ be a zeroth-order morphism from $\boldsymbol{\mathcal{A}}^{(0)}$ to $\boldsymbol{\mathcal{B}}^{(0)}$ and $\mathbf{g}^{(0)} = (\mathbf{d}, g^{(0)})$ a zeroth-order morphism from $\boldsymbol{\mathcal{B}}^{(0)}$ to $\boldsymbol{\mathcal{C}}^{(0)}$. Note that, $g^{(0)\sharp}_{\varphi} \circ f^{(0)}$ is an $S$-sorted mapping from $X$ to $(\T_{\Omega}(Z)_{\psi})_{\varphi}$. Thus, according to Proposition~\ref{PPropUniv}, there exists a unique $\Sigma$-homomorphism
$$
\left( g^{(0)\sharp}_{\varphi} \circ f^{(0)} \right)^{\sharp} \colon
\mathrm{T}_{\Sigma}(X)
\mor
(\mathrm{T}_{\Omega}(Z)_{\psi})_{\varphi}
$$
from $\mathbf{T}_{\Sigma}(X)$ to $(\mathbf{d} \circ \mathbf{c})^{\ast}_{\mathfrak{d}}(\mathbf{T}_{\Omega}(Z))$ such that
$(g^{(0)\sharp}_{\varphi} \circ f^{(0)})^{\sharp} \circ \eta^{X} = g^{(0)\sharp}_{\varphi} \circ f^{(0)}.$

However, considering the $S$-sorted mapping $f^{(0)}$ from $X$ to $\T_{\Lambda}(Y)_{\varphi}$, there exists a unique $\Sigma$-homomorphism 
$$
f^{(0)\sharp} \colon \mathrm{T}_{\Sigma}(X) \mor \mathrm{T}_{\Lambda}(Y)_{\varphi}
$$
from $\mathbf{T}_{\Sigma}(X)$ to $\mathbf{c}^{\ast}_{\mathfrak{d}}(\mathbf{T}_{\Lambda}(Y))$ such that $f^{(0)\sharp} \circ \eta^{X} = f^{(0)}$. Moreover, the same construction with $g^{(0)}$ leads to a $\Lambda$-homomorphism 
$$
g^{(0)\sharp} \colon \mathrm{T}_{\Lambda}(Y) \mor \mathrm{T}_{\Omega}(Z)_{\psi}
$$
from $\mathbf{T}_{\Lambda}(Y)$ to $\mathbf{d}_{\mathfrak{d}}^{\ast}(\mathbf{T}_{\Omega}(Z))$ such that $(g^{(0)\sharp} \circ \eta^{Y}) = g^{(0)}$. Thus, according to Proposition~\ref{PFunSig}, the mapping
$$
g^{(0)\sharp}_{\varphi} \colon \mathrm{T}_{\Lambda}(Y)_{\varphi} \mor (\mathrm{T}_{\Omega}(Z)_{\psi})_{\varphi}
$$
is a $\Sigma$-homomorphism from $\mathbf{c}_{\mathfrak{d}}^{\ast}(\mathbf{T}_{\Lambda}(Y))$ to $\mathbf{c}_{\mathfrak{d}}^{\ast}(\mathbf{d}_{\mathfrak{d}}^{\ast}(\mathbf{T}_{\Omega}(Z)))$.

All in all, since
$$
(g^{(0)\sharp}_{\varphi} \circ f^{(0)\sharp}) \circ \eta^{X}
=
g^{(0)\sharp}_{\varphi} \circ (f^{(0)\sharp} \circ \eta^{X})
=
g^{(0)\sharp}_{\varphi} \circ f^{(0)},
$$
it follows by uniqueness that
$
(g^{(0)\sharp}_{\varphi} \circ f^{(0)})^{\sharp}
$ 
is equal to 
$
g^{(0)\sharp}_{\varphi} \circ f^{(0)\sharp}.
$

This process is depicted in Figure~\ref{FCompExt}
\begin{figure}
$$
\xymatrix{
&&&
\mathrm{T}_{\Sigma}(X)
  \ar[d]_{f^{(0)\sharp}}
  \ar`l[dd]+/l4pc/`[dd]_{(g^{(0)\sharp}_{\varphi} \circ f^{(0)})^{\sharp}} [dd]
&
X
  \ar[l]_{\eta^{X}}
  \ar[ld]^{f^{(0)}}
  \ar@(d,r)[ldd]^{g^{(0)\sharp}_{\varphi} \circ f^{(0)}}
\\
Y
  \ar[r]^{\eta^{Y}}
  \ar[rd]_{g^{(0)}}
&
\mathrm{T}_{\Lambda}(Y)
  \ar[d]^{g^{(0)\sharp}}
&&
\mathrm{T}_{\Lambda}(Y)_{\varphi}
  \ar[d]_{g^{(0)\sharp}_{\varphi}}
\\
&
\mathrm{T}_{\Omega}(Z)_{\psi}
&&
(\mathrm{T}_{\Omega}(Z)_{\psi})_{\varphi}
}
$$
\caption{Construction of the $S$-sorted mapping $g^{(0)\sharp}_{\varphi} \circ f^{(0)\sharp}$.}
\label{FCompExt}
\end{figure}
\end{remark}

We now show that zeroth-order many-sorted rewriting systems and morphisms between them constitute a category.

\begin{proposition}\label{PRws0Cat}
The zeroth-order many-sorted rewriting systems together with zeroth-order morphisms constitute a category, denoted by $\mathsf{Rws}^{(0)}_{\mathfrak{d}}$.
\end{proposition}

\begin{proof}
That domains and codomains respect identities and compositions follows from the definitions of identity morphism and composition of morphisms introduced in Definitions~\ref{DIdRws0} and \ref{DCompRws0}. Thus, all that remains to be proven is that the identity morphism at a zeroth-order many-sorted rewriting system acts as a unit element and that the composition of morphisms of zeroth-order many-sorted rewriting systems is associative.

\textsf{Unit element.}

Let $\mathbf{f}^{(0)}=(\varphi, c, f^{(0)})$ be a morphism of zeroth-order many-sorted rewriting systems from $\boldsymbol{\mathcal{A}}^{(0)}=(S, \Sigma, X)$ to $\boldsymbol{\mathcal{B}}^{(0)}=(T, \Lambda, Y)$. We need to prove that
\begin{align*}
\mathbf{f}^{(0)}\circ\mathrm{id}^{\boldsymbol{\mathcal{A}}^{(0)}} &= \mathbf{f}^{(0)}
&&\mbox{ and }&
\mathrm{id}^{\boldsymbol{\mathcal{B}}^{(0)}}\circ\mathbf{f}^{(0)} &= \mathbf{f}^{(0)}.
\tag{Id}
\end{align*}

We will only prove the left-hand side of (Id). The right-hand side of (Id) is done similarly. According to Definition~\ref{DIdRws0}, the identity morphism at $\boldsymbol{\mathcal{A}}^{(0)}$ is given by $(\mathrm{id}^{S}, \eta^{\Sigma}, \eta^{X})$. Moreover, according to Definition~\ref{DCompRws0}, the composition zeroth-order morphisms is
$$
\mathbf{f}^{(0)}\circ\mathrm{id}^{\boldsymbol{\mathcal{A}}^{(0)}}
=
(\mathbf{c}, f^{(0)})\circ(\mathrm{id}^{\mathbf{\Sigma}}, \eta^{X})
=
\left(
\mathbf{c}\circ\mathrm{id}^{\mathbf{\Sigma}},
f^{(0)\sharp}_{\mathrm{id}^{S}} \circ \eta^{X}
\right).
$$

Note that, according Proposition~\ref{PSigdCat}, $\mathsf{Sig}_{\mathfrak{d}}$ is a category, therefore, $\mathbf{c}\circ\mathrm{id}^{\mathbf{\Sigma}}$ is equal to $\mathbf{c}$. Thus, all what remains to be proven is that
$$
f^{(0)\sharp}_{\mathrm{id}^{S}} \circ \eta^{X}
=
f^{(0)}.
$$

Note that the following chain of equalities holds
\allowdisplaybreaks
\begin{align*}
f^{(0)\sharp}_{\mathrm{id}^{S}} \circ \eta^{X}
&= 
\mathrm{Id}^{\mathsf{Alg}(\Sigma)} \left( f^{(0)\sharp} \right) \circ \eta^{X}
\tag{1}
\\
&=
f^{(0)\sharp} \circ \eta^{X}
\tag{2}
\\
&=
f^{(0)}.
\tag{3}
\end{align*}

The first equality follows from the fact that, by Proposition~\ref{PAlgFun}, $\mathrm{Alg}_{\mathfrak{d}}$ is a functor; the second equality unravels the definition of identity functor at $\mathsf{Alg}(\Sigma)$; finally, the last equality follows from Proposition~\ref{PPropUniv}.

Thus, Equation (Id) follows.

This proves that identities at zeroth-order many-sorted rewriting systems act as the unit element for composition of morphisms between zeroth-order many-sorted rewriting systems.

\textsf{Associativity.}

Let $\boldsymbol{\mathcal{A}}^{(0)}$, $\boldsymbol{\mathcal{B}}^{(0)}$, $\boldsymbol{\mathcal{C}}^{(0)}$ and $\boldsymbol{\mathcal{D}}^{(0)}$ be four zeroth-order many-sorted rewriting systems 
of the form
$$
\boldsymbol{\mathcal{A}}^{(0)}=(S, \Sigma, X), \,
\boldsymbol{\mathcal{B}}^{(0)}=(T, \Lambda, Y), \,
\boldsymbol{\mathcal{C}}^{(0)}=(U, \Omega, Z)
$$
$$
\mbox{and }
\boldsymbol{\mathcal{D}}^{(0)}=(V, \Gamma, W),
\mbox{ respectively},
$$
and let $\mathbf{f}^{(0)}$, $\mathbf{g}^{(0)}$ and $\mathbf{h}^{(0)}$ be zeroth-order morphisms of the form
$$
\mathbf{f}^{(0)} \colon \boldsymbol{\mathcal{A}}^{(0)} \mor \boldsymbol{\mathcal{B}}^{(0)}, \,
\mathbf{g}^{(0)} \colon \boldsymbol{\mathcal{B}}^{(0)} \mor \boldsymbol{\mathcal{C}}^{(0)}
\mbox{ and }
\mathbf{h}^{(0)} \colon \boldsymbol{\mathcal{C}}^{(0)} \mor \boldsymbol{\mathcal{D}}^{(0)}
$$
where $\mathbf{f}^{(0)}$, $\mathbf{g}^{(0)}$ and $\mathbf{h}^{(0)}$ stand for
$$
\mathbf{f}^{(0)}=(\varphi, c, f^{(0)}), \,
\mathbf{g}^{(0)}=(\psi, d, g^{(0)})
\mbox{ and }
\mathbf{h}^{(0)}=(\theta, e, h^{(0)}),
\mbox{ respectively.}
$$
We need to prove that 
\begin{equation}
\mathbf{h}^{(0)}\circ(\mathbf{g}^{(0)}\circ\mathbf{f}^{(0)})
=
(\mathbf{h}^{(0)}\circ\mathbf{g}^{(0)})\circ\mathbf{f}^{(0)}.
\tag{Assoc}
\end{equation}

Following Definition~\ref{DCompRws0}, the left-hand side composition of Equation (Assoc) is
\allowdisplaybreaks
\begin{align*}
\mathbf{h}^{(0)} \circ (\mathbf{g}^{(0)} \circ \mathbf{f}^{(0)})
&=
(\mathbf{e}, h^{(0)}) \circ \left( (\mathbf{d}, g^{(0)}) \circ (\mathbf{c}, f^{(0)}) \right)
\\
&=
(\mathbf{e}, h^{(0)}) \circ \left(
\mathbf{d}\circ\mathbf{c}, 
g^{(0)\sharp}_{\varphi}\circ f^{(0)}
\right)
\\
&= 
\left(
\mathbf{e}\circ(\mathbf{d}\circ\mathbf{c}),
h^{(0)\sharp}_{\psi \circ \varphi} \circ \left( g^{(0)\sharp}_{\varphi}\circ f^{(0)} \right)
\right).
\end{align*}
Similarly, the right-hand side composition of Equation~(Assoc) is
\allowdisplaybreaks
\begin{align*}
(\mathbf{h}^{(0)} \circ \mathbf{g}^{(0)}) \circ \mathbf{f}^{(0)}
&=
\left( (\mathbf{e}, h^{(0)}) \circ (\mathbf{d}, g^{(0)}) \right) \circ (\mathbf{c}, f^{(0)})
\\
&=
\left(
\mathbf{e}\circ\mathbf{d},
h^{(0)\sharp}_{\psi} \circ g^{(0)}
\right) \circ (\mathbf{c}, f^{(0)})
\\
&=
\left(
(\mathbf{e}\circ\mathbf{d})\circ\mathbf{c},
\left(h^{(0)\sharp}_{\psi} \circ g^{(0)}\right)^{\sharp}_{\varphi} \circ f^{(0)}
\right).
\end{align*}

Note that, according Proposition~\ref{PSigdCat}, $\mathsf{Sig}_{\mathfrak{d}}$ is a category, therefore, 
$$
\mathbf{e}\circ(\mathbf{d}\circ\mathbf{c})
=
(\mathbf{e}\circ\mathbf{d})\circ\mathbf{c}.
$$
Thus, all what remains to be proven is that
$$
h^{(0)\sharp}_{\psi \circ \varphi} \circ \left( g^{(0)\sharp}_{\varphi}\circ f^{(0)} \right)
=
\left(h^{(0)\sharp}_{\psi} \circ g^{(0)}\right)^{\sharp}_{\varphi} \circ f^{(0)}.
$$
Moreover, it suffices to show that
$$
h^{(0)\sharp}_{\psi \circ \varphi} \circ g^{(0)\sharp}_{\varphi}
=
\left(h^{(0)\sharp}_{\psi} \circ g^{(0)}\right)^{\sharp}_{\varphi}.
$$

The following chain of equalities holds
\allowdisplaybreaks
\begin{align*}
h^{(0)\sharp}_{\psi \circ \varphi} \circ g^{(0)\sharp}_{\varphi}
&=
(h^{(0)\sharp}_{\psi})_{\varphi} \circ g^{(0)\sharp}_{\varphi}
\tag{1}
\\
&=
(h^{(0)\sharp}_{\psi} \circ g^{(0)\sharp})_{\varphi}
\tag{2}
\\
&=
\left(h^{(0)\sharp}_{\psi} \circ g^{(0)} \right)^{\sharp}_{\varphi}.
\tag{3}
\end{align*}

The first equality follows from the fact that, by Proposition~\ref{PAlgFun}, $\mathrm{Alg}_{\mathfrak{d}}$ is a contravariant functor; the second equality follows from the fact that, by Proposition~\ref{PFunSig}, $\mathbf{c}^{\ast}_{\mathfrak{d}}$ is a covariant functor; finally, the last equality follows the fact that, according to Remark~\ref{RCompExt},
$
(h^{(0)\sharp}_{\psi} \circ g^{(0)\sharp})_{\varphi}
$
is equal to
$
(h^{(0)\sharp}_{\psi} \circ g^{(0)})^{\sharp}_{\varphi}.
$

Thus, Equation (Assoc) follows.

This proves that composition of morphisms between zeroth-order many-sorted rewriting systems is associative.

This proves that $\mathsf{Rws}^{(0)}_{\mathfrak{d}}$ is a category.
\end{proof}

We next associate, to every zeroth-order many-sorted rewriting system, a zeroth-order tower. Moreover, we define the concept of morphism between zeroth-order tower.

\begin{definition}\label{DTw0}
Given a zeroth-order many-sorted rewriting system $\boldsymbol{\mathcal{A}}^{(0)} = (S, \Sigma, X)$, the \emph{zeroth-order tower} associated to $\boldsymbol{\mathcal{A}}^{(0)}$ is $\mathbb{A}^{(0)}=(\boldsymbol{\mathcal{A}}^{(0)}, \mathbf{T}_{\Sigma}(X))$ where $\mathbf{T}_{\Sigma}(X)$ is the free $\Sigma$-algebra on $X$.

A \emph{morphism of zeroth-order towers} from $\mathbb{A}^{(0)} = (\boldsymbol{\mathcal{A}}^{(0)}, \T_{\Sigma}(X))$ to $\mathbb{B}^{(0)} = (\boldsymbol{\mathcal{B}}^{(0)}, \T_{\Lambda}(Y))$ is an ordered triple $(\mathbb{A}^{(0)}, \mathbf{f}^{\sharp(0)}, \mathbb{B}^{(0)})$, denoted by $\mathbf{f}^{\sharp(0)} \colon \mathbb{A}^{(0)} \mor \mathbb{B}^{(0)}$ for short, in which $\mathbf{f}^{\sharp (0)}$ is the ordered pair $(\mathbf{f}^{(0)}, f^{\sharp (0)})$ where $\mathbf{f}^{(0)} = (\varphi, c, f^{(0)})$ is a zeroth-order morphism from $\boldsymbol{\mathcal{A}}^{(0)}$ to $\boldsymbol{\mathcal{B}}^{(0)}$ and $f^{(0)\sharp}$ is the unique $\Sigma$-homomorphism from $\mathbf{T}_{\Sigma}(X)$ to $\mathbf{c}_{\mathfrak{d}}^{\ast}(\mathbf{T}_{\Lambda}(Y))$ introduced in Proposition~\ref{PPropUniv}. 
\end{definition}

We define below the notions of identity morphism at a zeroth-order tower and composite morphism of zeroth-order towers.

\begin{definition}\label{DIdTw0}
Let $\mathbb{A}^{(0)}=(\boldsymbol{\mathcal{A}}^{(0)}, \mathbf{T}_{\Sigma}(X))$ be the zeroth-order tower associated to the zeroth-order many-sorted rewriting system $\boldsymbol{\mathcal{A}}^{(0)}=(S, \Sigma, X)$. The ordered pair $(\mathrm{id}^{\boldsymbol{\mathcal{A}}^{(0)}}, \mathrm{id}^{\T_{\Sigma}(X)})$, denoted by $\mathrm{id}^{\mathbb{A}^{(0)}}$, is the identity morphism at $\mathbb{A}^{(0)}$. Recall that, $\mathrm{id}^{\boldsymbol{\mathcal{A}}^{(0)}}$ is the identity morphism at $\boldsymbol{\mathcal{A}}^{(0)}$ in $\mathsf{Rws}_{\mathfrak{d}}^{(0)}$ introduced at Definition~\ref{DIdRws0}. Moreover, the $\Sigma$-algebra $(\mathrm{id}^{\mathbf{\Sigma}})^{\ast}_{\mathfrak{d}}(\mathbf{T}_{\Sigma}(X))$ is $\mathbf{T}_{\Sigma}(X)$ since, by Proposition~\ref{PAlgFun}, $\mathrm{Alg}_{\mathfrak{d}}$ is a functor. Thus, the canonical extension $\left( \eta^{X} \right)^{\sharp}$ is $\mathrm{id}^{\T_{\Sigma}(X)}$.
\end{definition}

\begin{definition}\label{DCompTw0}
Let $\mathbf{f}^{(0)\sharp}=(\mathbf{f}^{(0)}, f^{(0)\sharp})$ be a morphism of zeroth-order towers from $\mathbb{A}^{(0)}$ to $\mathbb{B}^{(0)}$ and $\mathbf{g}^{(0)\sharp}=(\mathbf{g}^{(0)}, g^{(0)\sharp})$ a morphism of zeroth-order towers from $\mathbb{B}^{(0)}$ to $\mathbb{C}^{(0)}$ where $\mathbb{A}^{(0)}$, $\mathbb{B}^{(0)}$ and $\mathbb{C}^{(0)}$ stand for
$$
\mathbb{A}^{(0)}=(\boldsymbol{\mathcal{A}}^{(0)}, \mathbf{T}_{\Sigma}(X)), \,
\mathbb{B}^{(0)}=(\boldsymbol{\mathcal{B}}^{(0)}, \mathbf{T}_{\Lambda}(Y))
$$
$$
\mbox{ and }
\mathbb{C}^{(0)}=(\boldsymbol{\mathcal{C}}^{(0)}, \mathbf{T}_{\Omega}(Z)),
\mbox{ respectively.}
$$
The zeroth-order \emph{composition} tower morphism $\mathbf{g}^{(0)\sharp}\circ\mathbf{f}^{(0)\sharp}$, from $\mathbb{A}^{(0)}$ to $\mathbb{C}^{(0)}$, is
$$
\mathbf{g}^{(0)\sharp}\circ\mathbf{f}^{(0)\sharp}=
(\mathbf{g}^{(0)}, g^{(0)\sharp})\circ(\mathbf{f}^{(0)},f^{(0)\sharp})=
\left(
\mathbf{g}^{(0)}\circ\mathbf{f}^{(0)}, 
g^{(0)\sharp}_{\varphi} \circ f^{(0)\sharp}
\right).
$$
The composition $\mathbf{g}^{(0)}\circ\mathbf{f}^{(0)}$ is the morphism of $\mathsf{Rws}_{\mathfrak{d}}^{(0)}$ from $\boldsymbol{\mathcal{A}}^{(0)}$ to $\boldsymbol{\mathcal{C}}^{(0)}$ which, according to Definition~\ref{DCompRws0}, is
$$
\mathbf{g}^{(0)}\circ\mathbf{f}^{(0)}=
\left(
\mathbf{d}\circ\mathbf{c}, 
g^{(0)\sharp}_{\varphi} \circ f^{(0)}
\right).
$$
and, according to Remark~\ref{RCompExt},
$
g^{(0)\sharp}_{\varphi} \circ f^{(0)\sharp}.
$
its extension.
\end{definition}

We now show that zeroth-order towers and morphisms between them constitute a category.

\begin{proposition}\label{PTw0Cat}
The zeroth-order towers together with zeroth-order morphisms between zeroth-order towers constitute a category, denoted by $\mathsf{Tw}^{(0)}_{\mathfrak{d}}$.
\end{proposition}

\begin{proof}
That domains and codomains respect identities and compositions follows from the definitions of identity morphism and composition of morphisms introduced in Definitions~\ref{DIdTw0} and \ref{DCompTw0}. Thus, all that remains to be proven is that the identity morphism at a zeroth-order tower acts as a unit element and that the composition of morphisms of zeroth-order towers is associative.

{\sffamily Unit element.}

Let $\mathbf{f}^{(0)\sharp}=(\mathbf{f}^{(0)}, f^{(0)\sharp})$ be a zeroth-order morphism from $\mathbb{A}^{(0)}$ to $\mathbb{B}^{(0)}$. We need to prove that 
\begin{align*}
\mathbf{f}^{(0)\sharp} \circ \mathrm{id}^{\mathbb{A}^{(0)}} &= \mathbf{f}^{(0)\sharp}
&&\mbox{ and }&
\mathrm{id}^{\mathbb{B}^{(0)}} \circ \mathbf{f}^{(0)\sharp} &= \mathbf{f}^{(0)\sharp}.
\tag{Id}
\end{align*}

We will only prove the left-hand side of (Id). The right-hand side of (Id) is done similarly. Let us prove that $\mathbf{f}^{(0)\sharp} \circ \mathrm{id}^{\mathbb{A}^{(0)}} = \mathbf{f}^{(0)\sharp}$.

The following chain of equalities holds
\allowdisplaybreaks
\begin{align*}
\mathbf{f}^{(0)\sharp} \circ \mathrm{id}^{\mathbb{A}^{(0)}}
&=
\left( \mathbf{f}^{(0)}, f^{(0)\sharp} \right) \circ \left( \mathrm{id}^{\boldsymbol{\mathcal{A}}^{(0)}}, \mathrm{id}^{\T_{\Sigma}(X)} \right)
\tag{1}
\\
&=
\left( \mathbf{f}^{(0)}, f^{(0)\sharp} \right) \circ \left( \mathrm{id}^{\boldsymbol{\mathcal{A}}^{(0)}}, \left( \eta^{X} \right)^{\sharp} \right)
\tag{2}
\\
&=
\left( \mathbf{f}^{(0)} \circ \mathrm{id}^{\boldsymbol{\mathcal{A}}^{(0)}}, f^{(0)\sharp}_{\mathrm{id}^{S}} \circ \left( \eta^{X} \right)^{\sharp} \right)
\tag{3}
\\
&=
\left( \mathbf{f}^{(0)}, \left( f^{(0)\sharp}_{\mathrm{id}^{S}} \circ \eta^{X}\right)^{\sharp}\right)
\tag{4}
\\
&=
(\mathbf{f}^{(0)}, f^{(0)\sharp})
\tag{5}
\\
&=
\mathbf{f}^{(0)\sharp}.
\tag{6}
\end{align*}

The first equality unravels the definitions of $\mathbf{f}^{(0)\sharp}$ and $\mathrm{id}^{\mathbb{A}^{(0)}}$, introduced in Definition~\ref{DIdTw0}; the second equality follows from the fact that $(\eta^{X})^{\sharp}=\mathrm{Id}^{\T_{\Sigma}(X)}$; the third equality unravels the definition of the composition of morphisms of zeroth-order towers, introduced in Definition~\ref{DCompTw0}; the fourth equality follows from the fact that, by Proposition~\ref{PRws0Cat}, $\mathsf{Rws}_{\mathfrak{d}}^{(0)}$ is category. Let us also recall that, as it was shown in Remark~\ref{RCompExt}, 
$$
f^{(0)\sharp}_{\mathrm{id}^{S}} \circ \left( \eta^{X} \right)^{\sharp}
=
\left(f^{(0)\sharp}_{\mathrm{id}^{S}} \circ \eta^{X}\right)^{\sharp}.
$$
The fifth equality follows from the fact that, according to Proposition~\ref{PAlgFun}, $\mathrm{Alg}_{\mathfrak{d}}$ is a functor, thus
$$
f^{(0)\sharp}_{\mathrm{id}^{S}} \circ \eta^{X}
=
f^{(0)\sharp} \circ \eta^{X}
=
f^{(0)};
$$
finally, the last equality recovers the definition of $\mathbf{f}^{(0)\sharp}$.

Thus, Equation (Id) follows.

This proves that identities at zeroth-order towers are the unit element for composition of morphisms between zeroth-order towers.

{\sffamily Associativity.}

Let $\mathbb{A}^{(0)}$, $\mathbb{B}^{(0)}$, $\mathbb{C}^{(0)}$ and $\mathbb{D}^{(0)}$ be four zeroth-order towers 
and let $\mathbf{f}^{(0)\sharp}=(\mathbf{f}^{(0)}, f^{(0)\sharp})$, $\mathbf{g}^{(0)\sharp}=(\mathbf{g}^{(0)},g^{(0)\sharp})$ and $\mathbf{h}^{(0)\sharp}=(\mathbf{h}^{(0)},h^{(0)\sharp})$ be morphisms of the form
$$
\mathbf{f}^{(0)\sharp} \colon \mathbb{A}^{(0)} \mor \mathbb{B}^{(0)}, \,
\mathbf{g}^{(0)\sharp} \colon \mathbb{B}^{(0)} \mor \mathbb{C}^{(0)}
\mbox{ and }
\mathbf{h}^{(0)\sharp} \colon \mathbb{C}^{(0)} \mor \mathbb{D}^{(0)}
$$
where $\mathbf{f}^{(0)}$, $\mathbf{g}^{(0)}$ and $\mathbf{h}^{(0)}$ stand for
$$
\mathbf{f}^{(0)}=(\varphi, c, f^{(0)}), \,
\mathbf{g}^{(0)}=(\psi, d, g^{(0)})
\mbox{ and }
\mathbf{h}^{(0)}=(\theta, e, h^{(0)}),
\mbox{ respectively.}
$$
We need to prove that
\begin{equation}
\mathbf{h}^{(0)\sharp} \circ \left(\mathbf{g}^{(0)\sharp} \circ \mathbf{f}^{(0)\sharp}\right)
=
\left(\mathbf{h}^{(0)\sharp} \circ \mathbf{g}^{(0)\sharp}\right) \circ \mathbf{f}^{(0)\sharp}.
\tag{Assoc}
\end{equation}

The following chain of equalities holds
\allowdisplaybreaks
\begin{align*}
&
\mathbf{h}^{(0)\sharp} \circ \left(\mathbf{g}^{(0)\sharp} \circ \mathbf{f}^{(0)\sharp}\right)
\\
&=
\left( \mathbf{h}^{(0)},h^{(0)\sharp} \right) \circ \left(\left( \mathbf{g}^{(0)}, g^{(0)\sharp} \right) \circ \left( \mathbf{f}^{(0)}, f^{(0)\sharp} \right) \right)
\tag{1}
\\
&=
\left(
\mathbf{h}^{(0)} \circ \left(\mathbf{g}^{(0)} \circ \mathbf{f}^{(0)}\right),
h^{(0)\sharp}_{\psi \circ \varphi} \circ \left(g^{(0)\sharp}_{\varphi} \circ f^{(0)\sharp}\right)
\right)
\tag{2}
\\
&=
\left(
\mathbf{h}^{(0)} \circ \left(\mathbf{g}^{(0)} \circ \mathbf{f}^{(0)}\right),
h^{(0)\sharp}_{\psi \circ \varphi} \circ \left(g^{(0)\sharp}_{\varphi} \circ f^{(0)}\right)^{\sharp}
\right)
\tag{3}
\\
&=
\left(
\mathbf{h}^{(0)} \circ \left(\mathbf{g}^{(0)} \circ \mathbf{f}^{(0)}\right),
\left(h^{(0)\sharp}_{\psi \circ \varphi} \circ \left( g^{(0)\sharp}_{\varphi} \circ f^{(0)}\right)\right)^{\sharp}
\right)
\tag{4}
\\
&=
\left(
\mathbf{h}^{(0)} \circ \left(\mathbf{g}^{(0)} \circ \mathbf{f}^{(0)}\right),
\left(\left(h^{(0)\sharp}_{\psi} \circ g^{(0)}\right)^{\sharp}_{\varphi} \circ f^{(0)}\right)^{\sharp}
\right)
\tag{5}
\\
&=
\left(
\left(\mathbf{h}^{(0)} \circ \mathbf{g}^{(0)}\right) \circ \mathbf{f}^{(0)},
\left(\left(h^{(0)\sharp}_{\psi} \circ g^{(0)}\right)^{\sharp}_{\varphi} \circ f^{(0)}\right)^{\sharp}
\right)
\tag{6}
\\
&=
\left(
\left(\mathbf{h}^{(0)} \circ \mathbf{g}^{(0)}\right) \circ \mathbf{f}^{(0)},
\left(h^{(0)\sharp}_{\psi} \circ g^{(0)}\right)^{\sharp}_{\varphi} \circ f^{(0)\sharp}
\right)
\tag{7}
\\
&=
\left(
\left(\mathbf{h}^{(0)} \circ \mathbf{g}^{(0)}\right) \circ \mathbf{f}^{(0)},
\left(h^{(0)\sharp}_{\psi} \circ g^{(0)\sharp}\right)_{\varphi} \circ f^{(0)\sharp}
\right)
\tag{8}
\\
&=
\left( \left( \mathbf{h}^{(0)},h^{(0)\sharp} \right) \circ \left( \mathbf{g}^{(0)}, g^{(0)\sharp} \right) \right) \circ \left( \mathbf{f}^{(0)}, f^{(0)\sharp} \right)
\tag{9}
\\
&=
\left(\mathbf{h}^{(0)\sharp} \circ \mathbf{g}^{(0)\sharp}\right) \circ \mathbf{f}^{(0)\sharp}.
\tag{10}
\end{align*}

The first equality unravels the definition of the morphisms $\mathbf{f}^{(0)\sharp}$, $\mathbf{g}^{(0)\sharp}$ and $\mathbf{h}^{(0)\sharp}$; the second equality unravels the definition of composition of morphisms between zeroth-order towers, introduced in Definition~\ref{DCompTw0}; the third and fourth equality follows from the explicit definition of the canonical extension shown in Remark~\ref{RCompExt}; the fifth equality follows from the fact that, as it was shown in Proposition~\ref{PRws0Cat},
$$
h^{(0)\sharp}_{\psi \circ \varphi} \circ \left( g^{(0)\sharp}_{\varphi} \circ f^{(0)}\right)
=
\left(h^{(0)\sharp}_{\psi} \circ g^{(0)}\right)^{\sharp}_{\varphi} \circ f^{(0)};
$$
the sixth equality follows from the fact that, by Proposition~\ref{PRws0Cat}, $\mathsf{Rws}_{\mathfrak{d}}^{(0)}$ is a category; the seventh and eighth equality follows from the explicit definition of the canonical extension shown in Remark~\ref{RCompExt}; the ninth equality recovers the definition of composition of morphisms between zeroth-order towers; finally, the last equality recovers the definition of the morphisms $\mathbf{f}^{(0)\sharp}$, $\mathbf{g}^{(0)\sharp}$ and $\mathbf{h}^{(0)\sharp}$.

Thus, Equation (Assoc) follows.

This proves that composition of morphisms between zeroth-order towers is associative.

This shows that $\mathsf{Tw}_{\mathfrak{d}}^{(0)}$ is a category.
\end{proof}

We will prove below that the categories just defined of zeroth-order many-sorted rewriting systems and zeroth-order towers are isomorphic. In particular, we define the assignments $V^{(0)}$ and $U^{(0)}$, we prove that they are functors between the respective categories and that they are mutually inverse functors.

\begin{definition}\label{DV(0)}
We let $V^{(0)}$ stand for the assignment from $\mathsf{Rws}_{\mathfrak{d}}^{(0)}$ to $\mathsf{Tw}_{\mathfrak{d}}^{(0)}$ defined as follows:
\begin{enumerate}
\item
for every zeroth-order many-sorted rewriting system $\boldsymbol{\mathcal{A}}^{(0)}=(S, \Sigma, X)$, $V^{(0)}(\boldsymbol{\mathcal{A}}^{(0)})$ is the associated zeroth-order tower $\mathbb{A}^{(0)} = (\boldsymbol{\mathcal{A}}^{(0)}, \mathbf{T}_{\Sigma}(X))$, and
\item
for every morphism $\mathbf{f}^{(0)} = (\mathbf{c}, f^{(0)}) \colon \boldsymbol{\mathcal{A}}^{(0)} \mor \boldsymbol{\mathcal{B}}^{(0)}$, $V^{(0)}(\mathbf{f}^{(0)})$ is the associated morphism $\mathbf{f}^{(0)\sharp}=(\mathbf{f}^{(0)}, f^{(0)\sharp}) \colon \mathbb{A}^{(0)} \mor \mathbb{B}^{(0)}$.
\end{enumerate}
\end{definition}

\begin{proposition}\label{PV(0)Fun}
The assignment $V^{(0)}$ from $\mathsf{Rws}_{\mathfrak{d}}^{(0)}$ to $\mathsf{Tw}_{\mathfrak{d}}^{(0)}$ is a covariant functor.
\end{proposition}

\begin{proof}
That $V^{(0)}$ maps objects and morphisms of $\mathsf{Rws}_{\mathfrak{d}}^{(0)}$ to objects and morphisms of $\mathsf{Tw}_{\mathfrak{d}}^{(0)}$ follows from the definition of the assignment. Therefore, all that remains to be proven is that $V^{(0)}$ preserves identities and compositions.

{\sffamily Preservation of identities.}

Let $\boldsymbol{\mathcal{A}}^{(0)} = (S, \Sigma, X)$ be an object in $\mathsf{Rws}_{\mathfrak{d}}^{(0)}$. We need to prove that
$$
V^{(0)} \left(\mathrm{id}^{\boldsymbol{\mathcal{A}}^{(0)}}\right)
=
\mathrm{id}^{V^{(0)}(\boldsymbol{\mathcal{A}}^{(0)})}.
$$

The following chain of equalities holds
\allowdisplaybreaks
\begin{align*}
V^{(0)} \left(\mathrm{id}^{\boldsymbol{\mathcal{A}}^{(0)}}\right)
&=
\left( \mathrm{id}^{\boldsymbol{\mathcal{A}}^{(0)}}, \left(\eta^{X}\right)^{\sharp}\right)
\tag{1}
\\
&=
\left( \mathrm{id}^{\boldsymbol{\mathcal{A}}^{(0)}}, \mathrm{id}^{\T_{\Sigma}(X)}\right)
\tag{2}
\\
&=
\mathrm{id}^{(\boldsymbol{\mathcal{A}}^{0}, \mathbf{T}_{\Sigma}(X))}
\tag{3}
\\
&=
\mathrm{id}^{V^{(0)}(\boldsymbol{\mathcal{A}}^{0})}.
\tag{4}
\end{align*}

The first equality unravels the definition of the assignment $V^{(0)}$ on morphisms, introduced in Definition~\ref{DV(0)}. Recall that, by Definition~\ref{DIdRws0}, the identity morphism at $\boldsymbol{\mathcal{A}}^{(0)}$ is
$$
\mathrm{id}^{\boldsymbol{\mathcal{A}}^{(0)}}=(\mathrm{id}^{S}, \eta^{\Sigma}, \eta^{X}).
$$
The second equality follows from the fact that $\left(\eta^{X}\right)^{\sharp} = \mathrm{id}^{\T_{\Sigma}(X)}$; the third equality recovers the definition of the identity morphism at the zeroth-order tower $(\boldsymbol{\mathcal{A}}^{(0)}, \mathbf{T}_{\Sigma}(X))$, introduced in Definition~\ref{DIdTw0}; finally, the last equality recovers the definition of the assignment $V^{(0)}$ on objects.

This proves that $V^{(0)}$ preserves identities.

{\sffamily Preservation of compositions.}

Let $\mathbf{f}^{(0)}=(\varphi, c, f^{(0)}) \colon \boldsymbol{\mathcal{A}}^{(0)} \mor \boldsymbol{\mathcal{B}}^{(0)}$ and $\mathbf{g}^{(0)}=(\psi, d, g^{(0)}) \colon \boldsymbol{\mathcal{B}}^{(0)} \mor \boldsymbol{\mathcal{C}}^{(0)}$ be morphisms in $\mathsf{Rws}_{\mathfrak{d}}^{(0)}$. We need to prove that
$$
V^{(0)} \left(\mathbf{g}^{(0)}\circ\mathbf{f}^{(0)}\right)
=
V^{(0)} \left(\mathbf{g}^{(0)}\right) \circ V^{(0)} \left(\mathbf{f}^{(0)}\right).
$$

The following chain of equalities holds
\begin{align*}
V^{(0)} \left(\mathbf{g}^{(0)}\circ\mathbf{f}^{(0)}\right)
&=
\left(\mathbf{g}^{(0)}\circ\mathbf{f}^{(0)},
\left(g^{(0)\sharp}_{\varphi} \circ f^{(0)}\right)^{\sharp}
\right)
\tag{1}
\\
&=
\left(
\mathbf{g}^{(0)} \circ  \mathbf{f}^{(0)},
g^{(0)\sharp}_{\varphi} \circ f^{(0)\sharp}
\right)
\tag{2}
\\
&=
\left(\mathbf{g}^{(0)}, g^{(0)\sharp}\right) \circ \left(\mathbf{f}^{(0)}, f^{(0)\sharp}\right)
\tag{3}
\\
&=
V^{(0)} \left(\mathbf{g}^{(0)}\right) \circ V^{(0)} \left(\mathbf{f}^{(0)}\right).
\tag{4}
\end{align*}

The first equality unravels the definition of the assignment $V^{(0)}$ on morphisms, introduced in Definition~\ref{DV(0)}. Recall that, by Definition~\ref{DCompRws0}, the comsposite morphism $\mathbf{g}^{(0)} \circ \mathbf{f}^{(0)}$ is
$$
\mathbf{g}^{(0)}\circ\mathbf{f}^{(0)}=
\left(
\mathbf{d}\circ\mathbf{c}, 
g^{(0)\sharp}_{\varphi} \circ f^{(0)}
\right).
$$
The second equality follows from the fact that, as it was shown in Remark~\ref{RCompExt}, 
$$
\left(g^{(0)\sharp}_{\varphi} \circ f^{(0)}\right)^{\sharp}
=
g^{(0)\sharp}_{\varphi} \circ f^{\sharp};
$$
the third equality recovers the definition of the composite morphisms between zeroth-order towers, introduced in Definition~\ref{DCompTw0}; finally, the last equality recovers the definition of the assignment $V^{(0)}$ on morphisms.

This proves that $V^{(0)}$ preserves compositions.

This completes the proof. 
\end{proof}

\begin{definition}\label{DU(0)}
We let $U^{(0)}$ stand for the assignment from $\mathsf{Tw}_{\mathfrak{d}}^{(0)}$ to $\mathsf{Rws}_{\mathfrak{d}}^{(0)}$ defined as follows:
\begin{enumerate}
\item
for every zeroth-order tower $\mathbb{A}^{(0)} = (\boldsymbol{\mathcal{A}}^{(0)}, \mathbf{T}_{\Sigma}(X))$, $U^{(0)}(\mathbb{A}^{(0)})$ is its underlying zeroth-order many-sorted rewriting system $\boldsymbol{\mathcal{A}}^{(0)}$, and
\item
for every morphism $\mathbf{f}^{(0)\sharp} = (\mathbf{f}^{(0)}, f^{(0)\sharp}) \colon \mathbb{A}^{(0)} \mor \mathbb{B}^{(0)}$, $U^{(0)}(\mathbf{f}^{(0)\sharp})$ is the underlying morphism $\mathbf{f}^{(0)} \colon \boldsymbol{\mathcal{A}}^{(0)} \mor \boldsymbol{\mathcal{B}}^{(0)}$.
\end{enumerate}
\end{definition}

\begin{proposition}\label{PU(0)Fun}
The assignment $U^{(0)}$ from $\mathsf{Tw}_{\mathfrak{d}}^{(0)}$ to $\mathsf{Rws}_{\mathfrak{d}}^{(0)}$ is a covariant functor.
\end{proposition}

\begin{proof}
That $U^{(0)}$ maps objects and morphisms of $\mathsf{Tw}_{\mathfrak{d}}^{(0)}$ to objects and morphisms of $\mathsf{Rws}_{\mathfrak{d}}^{(0)}$ follows from the definition of the assignment. Therefore, all that remains to be proven is that $U^{(0)}$ preserves identities and compositions.

{\sffamily Preservation of identities.}

Let $\mathbb{A}^{(0)} = (\boldsymbol{\mathcal{A}}^{(0)}, \mathbf{T}_{\Sigma}(X))$ be an object in $\mathsf{Tw}_{\mathfrak{d}}^{(0)}$ where $\boldsymbol{\mathcal{A}}^{(0)} = (S, \Sigma, X)$ is the underlying zeroth-order many-sorted rewriting system. We need to prove that
$$
U^{(0)} \left(\mathrm{id}^{\mathbb{A}^{(0)}}\right)
=
\mathrm{id}^{U^{(0)}(\mathbb{A}^{(0)})}.
$$

The following chain of equalities holds
\begin{align*}
U^{(0)} \left(\mathrm{id}^{\mathbb{A}^{(0)}}\right)
&=
U^{(0)} \left(\mathrm{id}^{\boldsymbol{\mathcal{A}}^{(0)}}, \mathrm{id}^{\T_{\Sigma}(X)}\right)
\\
&=
\mathrm{id}^{\boldsymbol{\mathcal{A}}^{(0)}}
\tag{2}
\\
&=
\mathrm{id}^{U^{(0)}(\mathbb{A}^{(0)})}.
\tag{3}
\end{align*}

The first equality unravels the definition of the identity morphism at $\mathbb{A}^{(0)}$;
the second equality unravels the definition of the assignment $U^{(0)}$ on morphisms; finally, the last equality recovers the definition of the assignment $U^{(0)}$ on objects.

This proves that $U^{(0)}$ preserves identities.

\textsf{Preservation of compositions.}

Let $\mathbf{f}^{(0)\sharp} = (\mathbf{f}^{(0)}, f^{(0)\sharp}) \colon \mathbb{A}^{(0)} \mor \mathbb{B}^{(0)}$ and $\mathbf{g}^{(0)\sharp} = (\mathbf{g}^{(0)}, g^{(0)\sharp}) \colon \mathbb{B}^{(0)} \mor \mathbb{C}^{(0)}$ be morphisms in $\mathsf{Tw}_{\mathfrak{d}}^{(0)}$. We need to prove that
$$
U^{(0)} \left(\mathbf{g}^{(0)} \circ \mathbf{f}^{(0)}\right)
=
U^{(0)} \left(\mathbf{g}^{(0)}\right) \circ U^{(0)} \left(\mathbf{f}^{(0)}\right).
$$

The following chain of equalities holds
\begin{align*}
U^{(0)} \left(\mathbf{g}^{(0)\sharp} \circ \mathbf{f}^{(0)\sharp}\right)
&=
U^{(0)} \left( \mathbf{g}^{(0)} \circ \mathbf{f}^{(0)}, g^{(0)\sharp}_{\varphi} \circ f^{(0)\sharp}\right)
\tag{1}
\\
&=
\mathbf{g}^{(0)} \circ \mathbf{f}^{(0)}
\tag{2}
\\
&=
U^{(0)} \left(\mathbf{g}^{(0)\sharp}\right) \circ U^{(0)} \left(\mathbf{f}^{(0)\sharp}\right).
\tag{3}
\end{align*}

The first equality unravels the definition of the composite of morphisms between zeroth-order towers, introduced in Definition~\ref{DCompTw0}; the second equality unravels the definition of the assignment $U^{(0)}$ on morphisms, introduced in Definition~\ref{DU(0)}; finally, the last equality recovers the definition of the assignment $U^{(0)}$ on morphisms.

This proves that $U^{(0)}$ preserves compositions.

This completes the proof.
\end{proof}

\begin{proposition}\label{PV(0)U(0)Comp}
$U^{(0)} \circ V^{(0)} = \mathrm{Id}^{\mathsf{Rws}_{\mathfrak{d}}^{(0)}}$ and $V^{(0)} \circ U^{(0)} = \mathrm{Id}^{\mathsf{Tw}_{\mathfrak{d}}^{(0)}}$
\end{proposition}

\begin{proof}
For every object $\boldsymbol{\mathcal{A}}^{(0)} = (S, \Sigma, X)$ in $\mathsf{Rws}_{\mathfrak{d}}^{(0)}$, the following chain of equalities holds
\begin{align*}
U^{(0)} \left(V^{(0)} \left(\boldsymbol{\mathcal{A}}^{(0)}\right)\right)
&=
U^{(0)} \left(\boldsymbol{\mathcal{A}}^{(0)}, \mathbf{T}_{\Sigma}(X)\right)
\tag{1}
\\
&=
\boldsymbol{\mathcal{A}}^{(0)}.
\tag{2}
\end{align*}

The first equality unravels the definition of the functor $V^{(0)}$ on objects, introduced in Definition~\ref{DV(0)}; the second equality unravels the definition of the functor $U^{(0)}$ on objects, introduced in Definition~\ref{DU(0)}.

For every morphism $\mathbf{f}^{(0)} \colon \boldsymbol{\mathcal{A}}^{(0)} \mor \boldsymbol{\mathcal{B}}^{(0)}$ in $\mathsf{Rws}_{\mathfrak{d}}^{(0)}$, the following chain of equalities holds
\begin{align*}
U^{(0)} \left(V^{(0)} \left(\mathbf{f}^{(0)}\right)\right)
&=
U^{(0)} \left(\mathbf{f}^{(0)}, f^{(0)\sharp}\right)
\tag{1}
\\
&=
\mathbf{f}^{(0)}.
\tag{2}
\end{align*}

The first equality unravels the definition of the functor $V^{(0)}$ on morphisms, introduced in Definition~\ref{DV(0)}; the second equality unravels the definition of the functor $U^{(0)}$ on morphisms, introduced in Definition~\ref{DU(0)}.

For every object $\mathbb{A}^{(0)} = (\boldsymbol{\mathcal{A}}^{(0)}, \mathbf{T}_{\Sigma}(X))$ in $\mathsf{Tw}_{\mathfrak{d}}^{(0)}$, the following chain of equalities holds
\begin{align*}
V^{(0)} \left(U^{(0)} \left(\mathbb{A}^{(0)}\right)\right)
&=
V^{(0)} \left(\boldsymbol{\mathcal{A}}^{(0)}\right)
\tag{1}
\\
&=
(\boldsymbol{\mathcal{A}}^{(0)}, \mathbf{T}_{\Sigma}(X))
\tag{2}
\\
&=
\mathbb{A}^{(0)}.
\tag{3}
\end{align*}

The first equality unravels the definition of the functor $U^{(0)}$ on objects, introduced in Definition~\ref{DU(0)}; the second equality unravels the definition of the functor $V^{(0)}$ on objects, introduced in Definition~\ref{DV(0)}; finally, the last equality recovers the definition $\mathbb{A}^{(0)}$.

For every morphism $\mathbf{f}^{(0)\sharp} = (\mathbf{f}^{(0)}, f^{(0)\sharp}) \colon \mathbb{A}^{(0)} \mor \mathbb{B}^{(0)}$ in $\mathsf{Tw}_{\mathfrak{d}}^{(0)}$ with $\mathbf{f}^{(0)}$, the following chain of equalities holds
\begin{align*}
V^{(0)} \left(U^{(0)} \left(\mathbf{f}^{(0)\sharp}\right)\right)
&=
V^{(0)} \left(\mathbf{f}^{(0)}\right)
\tag{1}
\\
&=
\left(\mathbf{f}^{(0)}, f^{(0)\sharp}\right)
\tag{2}
\\
&=
\mathbf{f}^{(0)}.
\tag{3}
\end{align*}

The first equality unravels the definition of the functor $U^{(0)}$ on morphisms, introduced in Definition~\ref{DU(0)}; the second equality unravels the definition of the functor $V^{(0)}$ on morphisms, introduced in Definition~\ref{DV(0)}; finally, the last equality recovers the definition $\mathbf{f}^{(0)\sharp}$.

This completes the proof.
\end{proof}

\begin{corollary}\label{CRws0Tw0Equiv}
The categories $\mathsf{Rws}_{\mathfrak{d}}^{(0)}$ and $\mathsf{Tw}_{\mathfrak{d}}^{(0)}$ are isomorphic.
\end{corollary} 			
\chapter{Many-sorted partial algebras}\label{S0C}

In this section we gather together the basic facts about many-sorted partial algebras that we will need afterwards. Specifically, we define the notion of many-sorted partial algebra, we introduce the homomorphisms between many-sorted partial algebras and several classes of them, as e.g., the closed and the full homomorphisms. Moreover, we define, in connection with the homomorphisms and their properties, several classes of subobjects and quotient objects of the many-sorted partial algebras. We also consider universal solutions, as e.g., free completions, and, finally, the concept of equation and the varieties of many-sorted partial algebras. 

\section{Many-sorted partial algebras}

\begin{definition} A \emph{structure of partial} $\Sigma$-\emph{algebra on} an $S$-sorted set $A$ is a family $(F_{\mathbf{s},s})_{(\mathbf{s},s)\in S^{\star}\times S}$, denoted by $F$, where, for $(\mathbf{s},s)\in S^{\star}\times S$, $F_{\mathbf{s},s}$ is a mapping from $\Sigma_{\mathbf{s},s}$ to $\mathrm{Hom}_{\mathrm{p}}(A_{\mathbf{s}},A_{s})$. If $\mathrm{Dom}(F_{\mathbf{s},s}(\sigma)) = A_{\mathbf{s}}$, then $F_{\mathbf{s},s}(\sigma)$ is a \emph{total operation} from $A_{\mathbf{s}}$ to $A_{s}$. If  $\mathrm{Dom}(F_{\mathbf{s},s}(\sigma)) = \varnothing$, then $F_{\mathbf{s},s}(\sigma)$ is called a \emph{discrete operation} or an \emph{empty (partial) operation} from $A_{\mathbf{s}}$ to $A_{s}$. If $(\mathbf{s},s) = (\lambda,s)$, $\sigma\in \Sigma_{\lambda,s}$, and $\mathrm{Dom}(F_{\lambda,s}(\sigma))\neq\varnothing$, then $F_{\lambda,s}(\sigma)$ picks out an element of $A_{s}$. Therefore a partial operation $F_{\lambda,s}(\sigma)$ from $A_{\lambda}$ to $A_{s}$, also called a \emph{partial constant} $F_{\lambda,s}(\sigma)$ from $A_{\lambda}$ to $A_{s}$, is either discrete or distinguishes exactly an element of $A_{s}$.  A \emph{many-sorted partial} $\Sigma$-\emph{algebra} (or, to abbreviate, \emph{partial $\Sigma$-\emph{algebra}}) is a pair $(A,F)$, abbreviated to $\mathbf{A}$, where $A$ is an $S$-sorted set and $F$ a structure of partial $\Sigma$-algebra on $A$. For a pair $(\mathbf{s},s)\in S^{\star}\times S$ and a formal operation $\sigma\in \Sigma_{\mathbf{s},s}$, in order to simplify the notation, the partial operation $F_{\mathbf{s},s}(\sigma)$ from $A_{\mathbf{s}}$ to $A_{s}$ will be written as $F_{\sigma}$. In some cases, to avoid mistakes, we will denote by $F^{\mathbf{A}}$ the structure of partial $\Sigma$-algebra on $A$, and for $(\mathbf{s},s)\in S^{\star}\times S$ and $\Sigma\in \Sigma_{\mathbf{s},s}$, by $F^{\mathbf{A}}_{\sigma}$, or simply by $\sigma^{\mathbf{A}}$, the corresponding partial operation.

A $\Sigma$-\emph{homomorphism} (or, to abbreviate, \emph{homomorphism}) from $\mathbf{A}$ to $\mathbf{B}$, where $\mathbf{B} = (B,G)$, is a triple $(\mathbf{A},f,\mathbf{B})$, abbreviated to $f\colon \mathbf{A}\mor \mathbf{B}$, where $f$ is an $S$-sorted mapping from $A$ to $B$ such that, for every $(\mathbf{s},s)\in S^{\star}\times S$, every $\sigma\in \Sigma_{\mathbf{s},s}$, and every $(a_{j})_{j\in \bb{\mathbf{s}}}\in A_{\mathbf{s}}$, if $(a_{j})_{j\in \bb{\mathbf{s}}}\in\mathrm{Dom}(\sigma^{\mathbf{A}})$, then $(f_{s_{j}}(a_{j}))_{j\in \bb{\mathbf{s}}}\in\mathrm{Dom}(\sigma^{\mathbf{B}})$ and
$
f_{s}(\sigma^{\mathbf{A}}((a_{j})_{j\in \bb{\mathbf{s}}})) = \sigma^{\mathbf{B}}((f_{s_{j}}(a_{j})_{j\in \bb{\mathbf{s}}})).
$
For a partial $\Sigma$-algebra $\mathbf{A}$ the triple $(\mathbf{A},\mathrm{id}^{A},\mathbf{A})$, where $\mathrm{id}^{A}$ is the identity at $A$, which is a homomorphism from $\mathbf{A}$ to $\mathbf{A}$, will be denoted by $\mathrm{id}^{\mathbf{A}}$ and will be called the \emph{identity} homomorphism at $\mathbf{A}$.
We will denote by $\mathsf{PAlg}(\Sigma)$ the category of partial $\Sigma$-algebras and homomorphisms  and by $\mathrm{PAlg}(\Sigma)$ the set of objects  of $\mathsf{PAlg}(\Sigma)$.
\end{definition}

\begin{remark}\label{RDiscAdj}
The forgetful functor $G_{\mathsf{PAlg}(\Sigma)}$ from the category $\mathsf{PAlg}(\Sigma)$ to the category $\mathsf{Set}^{S}$ has a left adjoint, denoted by $\mathbf{D}_{\Sigma}$, which sends an $S$-sorted set $A$ to $\mathbf{D}_{\Sigma}(A)$, the \emph{discrete} many-sorted partial $\Sigma$-algebra on $A$, defined as follows: Its underlying $S$-sorted set is $A$, and, for every $(\mathbf{s},s)\in S^{\star}\times S$ and every $\sigma\in \Sigma_{\mathbf{s},s}$, $\mathrm{Dom}(\sigma^{\mathbf{D}_{\Sigma}(A)}) = \varnothing$. To verify that 
$\mathbf{D}_{\Sigma}\dashv G_{\mathsf{PAlg}(\Sigma)}$ it suffices to note that, for every many-sorted partial 
$\Sigma$-algebra $\mathbf{B}$, it happens that every $S$-sorted mapping $f$ from $X$ to the underlying $S$-sorted set of $\mathbf{B}$ is a homomorphism from $\mathbf{D}_{\Sigma}(X)$ to $\mathbf{B}$.  
\end{remark}

\begin{remark}
The category $\mathsf{Alg}(\Sigma)$ is a full subcategory of the category $\mathsf{PAlg}(\Sigma)$. We shall prove later that the canonical full embedding of $\mathsf{Alg}(\Sigma)$ into $\mathsf{PAlg}(\Sigma)$ has a left adjoint, the free completion of a partial $\Sigma$-algebra, which will play a key role in this paper. 
\end{remark}

Since in this work some of the fundamental results consists in the proof that certain homomorphisms between many-sorted partial algebras are isomorphisms, we next state a useful characterization of them. 
(which show some similarity with the homeomorphisms, i.e., bijective and bicontinuous applications, between topological spaces). Further characterisations of the isomorphisms using the notions of closed homomorphism and complete homomorphism will be included below.

\begin{proposition}\label{PPAlgIso}
Let $f$ be a homomorphism from $\mathbf{A}$ to $\mathbf{B}$. Then the following statements are equivalent:
\begin{enumerate}
\item $f$ is an isomorphism from $\mathbf{A}$ to $\mathbf{B}$, i.e., $f$ is a section and a retraction.
\item $f$ is bijective and $(\mathbf{B},f^{-1},\mathbf{A})$, abbreviated to $f^{-1}\colon \mathbf{B}\mor \mathbf{A}$, is also a homomorphism from $\mathbf{B}$ to $\mathbf{A}$.
\end{enumerate}
\end{proposition}

\begin{remark}
Let $\mathbf{A}$ be a partial $\Sigma$-algebra. Then $(\mathbf{D}_{\Sigma}(A),\mathrm{id}^{A},\mathbf{A})$ is a bijective homomorphism, but it is an isomorphism if, and only if, $\mathbf{A}$ is also discrete.
\end{remark}

To establish some of the results that follow, it is useful to have the concept of product of a family of many-sorted partial algebras.

\begin{definition}
Let $(\mathbf{A}^{i})_{i\in I}$ be a family of partial $\Sigma$-algebras. Then
\begin{enumerate}
\item The \emph{product}
      of $(\mathbf{A}^{i})_{i\in I}$, $\prod_{i\in I}\mathbf{A}^{i}$, is
      the partial $\Sigma$-algebra which has as $S$-sorted underlying set $\prod_{i\in I}A^{i}$, 
      and where, for every $(\mathbf{s},s)\in S^{\star}\times S$ and every $\sigma\in \Sigma_{\mathbf{s},s}$,
      the structural operation $\sigma^{\prod_{i\in I}\mathbf{A}^{i}}$ is defined as follows:
      \begin{enumerate}
      \item $\mathrm{Dom}(\sigma^{\prod_{i\in I}\mathbf{A}^{i}}) = \{(a_{j})_{j\in \bb{\mathbf{s}}}\in 
      (\prod_{i\in I}A^{i})_{\mathbf{s}}\mid \forall\,i\in I\, 
      ((a_{j}(i))_{j\in \bb{\mathbf{s}}}\in \mathrm{Dom}(\sigma^{\mathbf{A}^{i}}))\}$.
      \item
      $$
      \sigma^{\prod_{i\in I}\mathbf{A}^{i}}
      \nfunction
      {\mathrm{Dom}(\sigma^{\prod_{i\in I}\mathbf{A}^{i}})}
      {\prod_{i\in I}A^{i}_{s}}
      {(a_{j})_{j\in \bb{\mathbf{s}}}}
      {(\sigma^{\mathbf{A}^{i}}((a_{j}(i))_{j\in \bb{\mathbf{s}}})_{i\in I}}
      $$
      \end{enumerate}
\item For every $i\in I$, the $i$-th \emph{canonical projection}
      is the homomorphism from $\prod_{i\in I}\mathbf{A}^{i}$ to $\mathbf{A}^{i}$
      determined by the $S$-sorted mapping $\pr^{i}$ which, for every $s\in S$, is
      defined as follows
      $$
      \pr^{i}_{s}
      \nfunction
      {\prod_{i\in I}A^{i}_{s}}
      {A^{i}_{s}}
      {(a_{i})_{ i\in I}}
      {a_{i}}
      $$
\end{enumerate}
\end{definition}

\begin{proposition}\label{PPAlgProd}
Let $(\mathbf{A}^{i})_{i\in I}$ be a family of partial $\Sigma$-algebras.
Then the pair $(\prod_{i\in I}\mathbf{A}^{i},(\pr^{i})_{i\in I})$
is a product in $\mathsf{PAlg}(\Sigma)$.
\end{proposition}

\section{Subalgebras}
In what follows we define the notions of subalgebra, relative subalgebra, and weak subalgebra of a partial $\Sigma$-algebra $\mathbf{A}$ and the subalgebra generating operator for $\mathbf{A}$. Furthermore, we characterize the subalgebras, the relative subalgebras and the weak subalgebras by means of the closed homomorphisms, the full homomorphisms and the homomorphisms, respectively.

\begin{definition}\label{PPAlgSub}
Let $\mathbf{A}$ be a partial $\Sigma$-algebra and $X\subseteq A$. Given $(\mathbf{s},s)\in S^{\star}\times S$ and $\sigma\in\Sigma_{\mathbf{s},s}$, we will say that $X$ is \emph{closed under the partial operation} $\sigma^{\mathbf{A}}\colon A_{\mathbf{s}}\dmor A_{s}$ if, for every $(a_{j})_{j\in\bb{\mathbf{s}}}\in X_{\mathbf{s}}$, if $(a_{j})_{j\in\bb{\mathbf{s}}}\in\mathrm{Dom}(\sigma^{\mathbf{A}})$, then $\sigma^{\mathbf{A}}((a_{j})_{j\in\bb{\mathbf{s}}})\in X_{s}$. We will say that $X$ is a \emph{closed subset} of $\mathbf{A}$ if $X$ is closed under the partial operations of $\mathbf{A}$. We will denote by $\mathrm{Cl}(\mathbf{A})$ the set of all closed subsets of $\mathbf{A}$ (which is an algebraic closure system on $A$) and by  $\mathbf{Cl}(\mathbf{A})$ the algebraic lattice $(\mathrm{Cl}(\mathbf{A}),\subseteq)$.

We will say that a partial $\Sigma$-algebra $\mathbf{B} $ is a \emph{subalgebra} of $\mathbf{A}$ if $B$ is a closed subset of $\mathbf{A}$, for every $(\mathbf{s},s)\in S^{\star}\times S$ and every $\sigma\in\Sigma_{\mathbf{s},s}$, $\mathrm{Dom}(\sigma^{\mathbf{B}}) = \mathrm{Dom}(\sigma^{\mathbf{A}})\cap B_{\mathbf{s}}$ and, for every  $(b_{j})_{j\in\bb{\mathbf{s}}}\in \mathrm{Dom}(\sigma^{\mathbf{B}})$, $\sigma^{\mathbf{B}}((b_{j})_{j\in\bb{\mathbf{s}}}) = \sigma^{\mathbf{A}}((b_{j})_{j\in\bb{\mathbf{s}}})$. We will denote by $\mathrm{Sub}(\mathbf{A})$ the set of all subalgebras of $\mathbf{A}$ and by $\mathbf{Sub}(\mathbf{A})$ the ordered set $(\mathrm{Sub}(\mathbf{A}),\subseteq)$. Since $\mathbf{Cl}(\mathbf{A})$ and $\mathbf{Sub}(\mathbf{A})$ are isomorphic, we shall feel free to deal either with a closed subset of $\mathbf{A}$ or with the correlated subalgebra of $\mathbf{A}$, whichever is most convenient for the work at hand.
\end{definition}

\begin{remark}
Let $\mathbf{A}$ and $\mathbf{B}$ be partial $\Sigma$-algebras such that $B\subseteq A$. Then $\mathbf{B} $ is a subalgebra of $\mathbf{A}$ if and only if, for every $(\mathbf{s},s)\in S^{\star}\times S$ and every $\sigma\in\Sigma_{\mathbf{s},s}$, $\Gamma_{\!\sigma^{\mathbf{B}}} = \Gamma_{\!\sigma^{\mathbf{A}}}\cap (B_{\mathbf{s}}\times A_{s})$.
\end{remark}

\begin{remark}
Let $\mathbf{A}$ be a partial $\Sigma$-algebra. If, for every $s\in S$ and every $\sigma\in \Sigma_{\lambda,s}$, $\mathrm{Dom}(\sigma^{\mathbf{A}}) = \varnothing$, then $\varnothing^{S}\in \mathrm{Cl}(\mathbf{A})$. Let us note that, in such a case, $\boldsymbol{\varnothing}^{S}$, the subalgebra of $\mathbf{A}$ canonically associated to $\varnothing^{S}$, is such that, for every $(\mathbf{s},s)\in S^{\star}\times S$ and every $\sigma\in\Sigma_{\mathbf{s},s}$, the partial operation from $\varnothing^{S}_{\mathbf{s}}$ to $\varnothing^{S}_{s} = \varnothing$ associated to $\sigma$ is discrete.
\end{remark}

We next define the notion of closed homomorphism from a partial $\Sigma$-algebra to another which will allow us, among other things, to provide another characterization of the isomorphisms and of the subalgebras of a partial $\Sigma$-algebra.

\begin{definition}
Let $\mathbf{A}$ and $\mathbf{B}$ be partial $\Sigma$-algebras and $f$ a homomorphism from $\mathbf{A}$ to $\mathbf{B}$. We will say that $f$ is \emph{closed} if, for every $(\mathbf{s},s)\in S^{\star}\times S$, every $\sigma\in\Sigma_{\mathbf{s},s}$, and every $(a_{j})_{j\in\bb{\mathbf{s}}}\in A_{\mathbf{s}}$, if $(f_{s_{j}}(a_{j})_{j\in\bb{\mathbf{s}}})\in \mathrm{Dom}(\sigma^{\mathbf{B}})$, then $(a_{j})_{j\in\bb{\mathbf{s}}}\in \mathrm{Dom}(\sigma^{\mathbf{A}})$, i.e., $\mathrm{Dom}(\sigma^{\mathbf{A}}) \supseteq f_{\mathbf{s}}^{-1}[\mathrm{Dom}(\sigma^{\mathbf{B}})]$ (therefore $\mathrm{Dom}(\sigma^{\mathbf{A}}) = f_{\mathbf{s}}^{-1}[\mathrm{Dom}(\sigma^{\mathbf{B}})]$).
\end{definition}

\begin{remark}
A homomorphism is closed if and only if the domains of the partial operations in the source are exactly the inverse images of the domains of the corresponding operations in the target.
\end{remark}

\begin{remark}
Every homomorphism between $\Sigma$-algebras is a closed homomorphism.
\end{remark}

\begin{proposition}
Let $f$ be a homomorphism from $\mathbf{A}$ to $\mathbf{B}$. Then the following statements are equivalent:
\begin{enumerate}
\item $f$ is an isomorphism from $\mathbf{A}$ to $\mathbf{B}$, i.e., $f$ is a section and a retraction.
\item $f$ is bijective and closed, i.e., for every $(\mathbf{s},s)\in S^{\star}\times S$, every $\sigma\in \Sigma_{\mathbf{s},s}$, and every $(a_{j})_{j\in \bb{\mathbf{s}}}\in A_{\mathbf{s}}$, the following holds: $(a_{j})_{j\in \bb{\mathbf{s}}}\in \mathrm{Dom}(\sigma^{\mathbf{A}})$ if, and only if, $ (f_{s_{j}}(a_{j}))_{j\in \bb{\mathbf{s}}}\in \mathrm{Dom}(\sigma^{\mathbf{B}})$.
\end{enumerate}
\end{proposition}

\begin{proposition}
Let $\mathbf{A}$ and $\mathbf{B}$ be partial $\Sigma$-algebras, $f$ a homomorphism from $\mathbf{A}$ to $\mathbf{B}$, and $X$ and $Y$ closed subsets of $\mathbf{A}$ and $\mathbf{B}$, respectively. Then:
\begin{enumerate}
\item $f^{-1}[Y]$ is a closed subset of $\mathbf{A}$.
\item If $f$ is closed, then $f[X]$ is a closed subset of $\mathbf{B}$.
\end{enumerate}
\end{proposition}

\begin{remark}
For a homomorphism $f$ from a partial $\Sigma$-algebra $\mathbf{A}$ to another $\mathbf{B}$ and a subset $X$ of $A$, we have that $f[\mathrm{Sg}_{\mathbf{A}}(X)] \subseteq \mathrm{Sg}_{\mathbf{B}}(f[X])$ (the converse inclusion does not holds in general; however, if $f$ is closed, then $f[\mathrm{Sg}_{\mathbf{A}}(X)] = \mathrm{Sg}_{\mathbf{B}}(f[X])$.
\end{remark}

We now characterize the subalgebras of a partial $\Sigma$-algebra by means of the closed homomorphisms.

\begin{proposition}
Let $\mathbf{A}$ and $\mathbf{B}$ be partial $\Sigma$-algebras such that $B\subseteq A$. Then $\mathbf{B}$ is a subalgebra of $\mathbf{A}$ if and only if $\mathrm{in}^{\mathbf{B},\mathbf{A}} = (\mathbf{B},\mathrm{in}^{B,A},\mathbf{A})$ is a closed homomorphism from $\mathbf{B}$ to $\mathbf{A}$.
\end{proposition}

\begin{definition}
Let $\mathbf{A}$ be a partial $\Sigma$-algebra. Then we denote by $\mathrm{Sg}_{\mathbf{A}}$ the algebraic closure operator canonically associated to the algebraic closure system $\mathrm{Cl}(\mathbf{A})$ on $A$ and we call it the \emph{subalgebra  generating operator for} $\mathbf{A}$. Moreover, if $X\subseteq A$, then we call $\mathrm{Sg}_{\mathbf{A}}(X)$ the \emph{subalgebra of} $\mathbf{A}$ \emph{generated by} $X$, and if $X$ is such that $\mathrm{Sg}_{\mathbf{A}}(X) = A$, then we will say that $X$ is a \emph{generating subset of} $\mathbf{A}$. Besides, $\mathbf{Sg}_{\mathbf{A}}(X)$ denotes the partial  $\Sigma$-algebra determined by $\mathrm{Sg}_{\mathbf{A}}(X)$.
\end{definition}

\begin{remark}
Let $\mathbf{A}$ be a partial $\Sigma$-algebra. Then the algebraic closure operator $\mathrm{Sg}_{\mathbf{A}}$ is uniform, i.e., for every  $X$, $Y\subseteq A$, if $\mathrm{supp}_{S}(X) = \mathrm{supp}_{S}(Y)$, then we have $\mathrm{supp}_{S}(\mathrm{Sg}_{\mathbf{A}}(X)) = \mathrm{supp}_{S}(\mathrm{Sg}_{\mathbf{A}}(Y))$.
\end{remark}

We next recall the proof by algebraic induction for partial $\Sigma$-algebras.

\begin{proposition}[Principle of Proof by Algebraic Induction]\label{PPAII}
Let $\mathbf{A}$ be a partial $\Sigma$-algebra generated by $X$. Then to prove that a subset $Y$ of $A$ is equal to $A$ it suffices to show: (1) $X\subseteq Y$ (algebraic induction basis) and (2) $Y$ is a closed subset subalgebra of $\mathbf{A}$ (algebraic induction step).
\end{proposition}

We next state the principle of extension of identities. This principle, which is fundamental to elucidate the equality of two coterminal homomorphisms, wil be used on several occasions in this work.

\begin{proposition}
Let $f,g\colon \mathbf{A}\mor \mathbf{B}$ be homomorphisms between partial $\Sigma$-algebras and let $X$ be a subset of $\mathbf{A}$. If $f\!\!\upharpoonright_{X} = g\!\!\upharpoonright_{X}$, then $f\!\!\upharpoonright_{\mathbf{Sg}_{\mathbf{A}}(X)} = g\!\!\upharpoonright_{\mathbf{Sg}_{\mathbf{A}}(X)}$. In particular, if $X$ is a generating subset of $\mathbf{A}$ and $f\!\!\upharpoonright_{X} = g\!\!\upharpoonright_{X}$, then $f = g$.
\end{proposition}

\begin{definition}\label{PPAlgSubRel}
Let $\mathbf{A}$ be a partial $\Sigma$-algebra. We will say that a partial $\Sigma$-algebra $\mathbf{B} $ is a \emph{relative subalgebra} of $\mathbf{A}$ if $B\subseteq A$ and, for every $(\mathbf{s},s)\in S^{\star}\times S$ and every $\sigma\in\Sigma_{\mathbf{s},s}$, $\mathrm{Dom}(\sigma^{\mathbf{B}}) \subseteq \mathrm{Dom}(\sigma^{\mathbf{A}})$ and, for every  $(b_{j})_{j\in\bb{\mathbf{s}}}\in B_{\mathbf{s}}$, if $(b_{j})_{j\in\bb{\mathbf{s}}}\in\mathrm{Dom}(\sigma^{\mathbf{A}})$ and $\sigma^{\mathbf{A}}((b_{j})_{j\in\bb{\mathbf{s}}})\in B_{s}$, then $(b_{j})_{j\in\bb{\mathbf{s}}}\in\mathrm{Dom}(\sigma^{\mathbf{B}})$ and $\sigma^{\mathbf{B}}((b_{j})_{j\in\bb{\mathbf{s}}}) = \sigma^{\mathbf{A}}((b_{j})_{j\in\bb{\mathbf{s}}})$. We will denote by $\mathrm{Sub}_{\mathrm{r}}(\mathbf{A})$ the set of all relative subalgebras of $\mathbf{A}$.
\end{definition}

\begin{remark}
Let $\mathbf{A}$ and $\mathbf{B}$ be partial $\Sigma$-algebras such that $B\subseteq A$. Then $\mathbf{B} $ is a relative subalgebra of $\mathbf{A}$ if and only if, for every $(\mathbf{s},s)\in S^{\star}\times S$ and every $\sigma\in\Sigma_{\mathbf{s},s}$, $\Gamma_{\!\sigma^{\mathbf{B}}} = \Gamma_{\!\sigma^{\mathbf{A}}}\cap (B_{\mathbf{s}}\times B_{s})$.
\end{remark}

\begin{remark}\label{PPAlgLift}
Let $\mathbf{A} = (A,F)$ be a partial $\Sigma$-algebra and $X\subseteq A$. Then there exists a unique structure of partial $\Sigma$-algebra $\mathrm{L}^{\mathrm{in}^{X,A}}(F)$ on $X$, the \emph{optimal lift of} $F$ \emph{with respect to } $\mathrm{in}^{X,A}$, such that, for every partial $\Sigma$-algebra $\mathbf{C}$ and every $S$-sorted mapping $h$ from $C$ to $X$, if $\mathrm{in}^{X,A}\circ h$ is a homomorphism from $\mathbf{C}$ to $\mathbf{A}$, then $h$ is a homomorphism from $\mathbf{C}$ to $(X,\mathrm{L}^{\mathrm{in}^{X,A}}(F))$. In fact, it suffices to take, for every $(\mathbf{s},s)\in S^{\star}\times S$ and every $\sigma\in \Sigma_{\mathbf{s},s}$, as partial operation associated to $\sigma$ the partial mapping $\mathrm{L}^{\mathrm{in}^{X,A}}(F)_{\sigma}$ from $X_{\mathbf{s}}$ to $X_{s}$ defined as follows:
$$
\mathrm{Dom}(\mathrm{L}^{\mathrm{in}^{X,A}}(F)_{\sigma}) = \{(x_{j})_{j\in\bb{\mathbf{s}}}\in X_{\mathbf{s}}\mid (x_{j})_{j\in\bb{\mathbf{s}}}\in \mathrm{Dom}(\sigma^{\mathbf{A}}) \And \sigma^{\mathbf{A}}((x_{j})_{j\in\bb{\mathbf{s}}})\in X_{s}\}
$$
and, for every $(x_{j})_{j\in\bb{\mathbf{s}}}\in \mathrm{Dom}(\mathrm{L}^{\mathrm{in}^{X,A}}(F)_{\sigma})$, $\mathrm{L}^{\mathrm{in}^{X,A}}(F)_{\sigma}((x_{j})_{j\in\bb{\mathbf{s}}}) = \sigma^{\mathbf{A}}((x_{j})_{j\in\bb{\mathbf{s}}})$. This shows that every subset $X$ of $A$ is the underlying $S$-sorted set of exactly one relative subalgebra of $\mathbf{A}$, and that the structure of partial $\Sigma$-algebra on $X$ is uniquely determined by the structure of partial $\Sigma$-algebra on $A$ and the specification of the subset $X$ of $A$.

Let us note that if the partial $\Sigma$-algebra $\mathbf{A}$ is \emph{total}, then, for every $(\mathbf{s},s)\in S^{\star}\times S$ and every $\sigma\in \Sigma_{\mathbf{s},s}$, we have that
$$
\mathrm{Dom}(\mathrm{L}^{\mathrm{in}^{X,A}}(F)_{\sigma}) = \{(x_{j})_{j\in\bb{\mathbf{s}}}\in X_{\mathbf{s}}\mid \sigma^{\mathbf{A}}((x_{j})_{j\in\bb{\mathbf{s}}})\in X_{s}\} = 
X_{\mathbf{s}}\cap (\sigma^{\mathbf{A}})^{-1}[X_{s}].
$$
Therefore, since $X_{\mathbf{s}}\cap (\sigma^{\mathbf{A}})^{-1}[X_{s}]\subseteq X_{\mathbf{s}}$ (but not necessarily $X_{\mathbf{s}}\subseteq (\sigma^{\mathbf{A}})^{-1}[X_{s}]$), we have that $(X,\mathrm{L}^{\mathrm{in}^{X,A}}(F))$ will, in general, be a \emph{partial} $\Sigma$-algebra.
\end{remark}

We next define the notion of full homomorphism from a partial $\Sigma$-algebra to another which will allow us, among other things, to provide another characterization of the isomorphisms and of the relative subalgebras of a $\Sigma$-algebra.

\begin{definition}
Let $\mathbf{A}$ and $\mathbf{B}$ be partial $\Sigma$-algebras and $f$ a homomorphism from $\mathbf{A}$ to $\mathbf{B}$. We will say that $f$ is \emph{full} if, for every $(\mathbf{s},s)\in S^{\star}\times S$, every $\sigma\in\Sigma_{\mathbf{s},s}$, and every $(a_{j})_{j\in\bb{\mathbf{s}}}\in A_{\mathbf{s}}$, if $((f_{s_{j}}(a_{j}))_{j\in\bb{\mathbf{s}}})\in \mathrm{Dom}(\sigma^{\mathbf{B}})$ and $\sigma^{\mathbf{B}}((f_{s_{j}}(a_{j}))_{j\in\bb{\mathbf{s}}})\in f_{s}[A_{s}]$, then there exists an $(a'_{j})_{j\in\bb{\mathbf{s}}}\in \mathrm{Dom}(\sigma^{\mathbf{A}})$ such that $((f_{s_{j}}(a'_{j}))_{j\in\bb{\mathbf{s}}}) = ((f_{s_{j}}(a_{j}))_{j\in\bb{\mathbf{s}}})$.
\end{definition}

\begin{remark}\label{PPAlgFullIni}
Those homomorphisms that are full and injective are characterizable as those that are initial and injective, where a homomorphism $f$ from a partial $\Sigma$-algebra $\mathbf{A}$ to another $\mathbf{B}$ is said to be \emph{initial} if, for every partial $\Sigma$-algebra $\mathbf{C}$ and every S-sorted mapping $h$ from $C$ to $A$, if $f\circ h$ is a homomorphism from $\mathbf{C}$ to $\mathbf{B}$, then $h$ is a homomorphism from $\mathbf{C}$ to $\mathbf{A}$. 
\end{remark}

\begin{proposition}
Let $f$ be a homomorphism from $\mathbf{A}$ to $\mathbf{B}$. Then the following statements are equivalent:
\begin{enumerate}
\item $f$ is an isomorphism from $\mathbf{A}$ to $\mathbf{B}$, i.e., $f$ is a section and a retraction.
\item $f$ is bijective and full.
\end{enumerate}
\end{proposition}

\begin{proposition}
Let $\mathbf{A}$ and $\mathbf{B}$ be partial $\Sigma$-algebras such that $B\subseteq A$. Then $\mathbf{B}$ is a relative subalgebra of $\mathbf{A}$ if and only if $\mathrm{in}^{\mathbf{B},\mathbf{A}} = (\mathbf{B},\mathrm{in}^{B,A},\mathbf{A})$ is a full homomorphism from $\mathbf{B}$ to $\mathbf{A}$.
\end{proposition}

%
%

\begin{definition}\label{PPAlgSubWeak}
Let $\mathbf{A}$ be a partial $\Sigma$-algebra. We will say that a partial $\Sigma$-algebra $\mathbf{B} $ is a \emph{weak subalgebra} of $\mathbf{A}$ if $B\subseteq A$ and, for every $(\mathbf{s},s)\in S^{\star}\times S$, every $\sigma\in\Sigma_{\mathbf{s},s}$, and every  $(b_{j})_{j\in\bb{\mathbf{s}}}\in B_{\mathbf{s}}$, if $(b_{j})_{j\in\bb{\mathbf{s}}}\in\mathrm{Dom}(\sigma^{\mathbf{B}})$, then $(b_{j})_{j\in\bb{\mathbf{s}}}\in\mathrm{Dom}(\sigma^{\mathbf{A}})$ and $\sigma^{\mathbf{B}}((b_{j})_{j\in\bb{\mathbf{s}}}) = \sigma^{\mathbf{A}}((b_{j})_{j\in\bb{\mathbf{s}}})$. We will denote by $\mathrm{Sub}_{\mathrm{w}}(\mathbf{A})$ the set of all weak subalgebras of $\mathbf{A}$.
\end{definition}

\begin{remark}
Let $\mathbf{A}$ and $\mathbf{B}$ be partial $\Sigma$-algebras such that $B\subseteq A$. Then $\mathbf{B}$ is a weak subalgebra of $\mathbf{A}$ if and only if, for every $(\mathbf{s},s)\in S^{\star}\times S$ and every $\sigma\in\Sigma_{\mathbf{s},s}$, $\Gamma_{\!\sigma^{\mathbf{B}}} \subseteq \Gamma_{\!\sigma^{\mathbf{A}}}$.
\end{remark}

\begin{proposition}
Let $\mathbf{A}$ and $\mathbf{B}$ be partial $\Sigma$-algebras such that $B\subseteq A$. Then $\mathbf{B}$ is a weak subalgebra of $\mathbf{A}$ if and only if $\mathrm{in}^{\mathbf{B},\mathbf{A}} = (\mathbf{B},\mathrm{in}^{B,A},\mathbf{A})$ is a homomorphism from $\mathbf{B}$ to $\mathbf{A}$.
\end{proposition}

\begin{remark}
Since every closed homomorphism is a full homomorphism and every full homomorphism is a homomorphism, it follows that every subalgebra is a relative subalgebra and every relative subalgebra is a weak subalgebra.  Moreover, every relative subalgebra on a closed subset is a subalgebra.
\end{remark}

\begin{proposition}
Let $\mathbf{A}$ and $\mathbf{B}$ be partial $\Sigma$-algebras and $f\colon A\mor B$ a mapping from $A$ to $B$. Then:
\begin{enumerate}
\item If $f\colon \mathbf{A}\mor \mathbf{B}$ is a homomorphism, then $\Gamma\!_{f}$ is a closed subset of $\mathbf{A}\times \mathbf{B}$.
\item If $\mathbf{B}$ is a $\Sigma$-algebra and $\Gamma\!_{f}$ is a closed subset of $\mathbf{A}\times \mathbf{B}$, then $f$ is a homomorphism from $\mathbf{A}$ to $\mathbf{B}$.
\end{enumerate}
\end{proposition}

\begin{remark}
Notice the (partial) analogy of the just stated proposition and the closed graph theorem in point-set topology.
\end{remark}

\begin{lemma}
Let $f\colon \mathbf{A}\mor \mathbf{B}$ be a homomorphism, $X\subseteq A$, $g\colon X\mor B$ such that $\Gamma\!_{g}\subseteq \Gamma\!_{f}$, and $\mathbf{X}$ the relative subalgebra of $\mathbf{A}$ on $X$. Then:
\begin{enumerate}
\item $g\colon \mathbf{X}\mor \mathbf{B}$ is a homomorphism.
\item $\Gamma\!_{g}$ is a closed subset of $\mathbf{A}\times \mathbf{B}$ if and only if $X$ is a closed subset of $\mathbf{A}$.
\end{enumerate}
\end{lemma}

\begin{corollary}
Let $f\colon \mathbf{A}\mor \mathbf{B}$ be a homomorphism and $X\subseteq A$. Then the following statements are equivalent:
\begin{enumerate}
\item $\mathrm{Sg}_{\mathbf{A}}(X) = A$.
\item $\mathrm{Sg}_{\mathbf{A}\times \mathbf{B}}(\Gamma\!_{f\upharpoonright_{X}}) = \Gamma\!_{f}$.
\end{enumerate}
\end{corollary}

We state the following lemma because it is used, along with other results, in the proof of the generalized recursion theorem, which we will also call Schmidt's construction.

\begin{lemma}
Let $\mathbf{A}$, $\mathbf{B}$, $\mathbf{C}$, and $\mathbf{D}$ be partial $\Sigma$-algebras, $\mathbf{C}$  a weak subalgebra of $\mathbf{A}$ such that $\mathrm{Sg}_{\mathbf{A}}(C) = A$, $\mathbf{D}$ a relative subalgebra of $\mathbf{B}$, and $f$ a homomorphism from $\mathbf{C}$ to $\mathbf{D}$ which allows a homomorphic extension $\overline{f}$ from $\mathbf{A}$ to $\mathbf{B}$. Then:
\begin{enumerate}
\item $\mathrm{Sg}_{\mathbf{A}\times\mathbf{D}}(\Gamma\!_{f})$ is the underlying function of a
      homomorphism $f^{\flat}$ from
      $$\mathbf{Dom}(\mathrm{Sg}_{\mathbf{A}\times\mathbf{D}}(\Gamma\!_{f}))$$ to $\mathbf{D}$
      (where $\mathbf{Dom}(\mathrm{Sg}_{\mathbf{A}\times\mathbf{D}}(\Gamma\!_{f}))$ is the
      relative subalgebra of $\mathbf{A}$ whose underlying many-sorted set is
      $\mathrm{Dom}(\mathrm{Sg}_{\mathbf{A}\times\mathbf{D}}(\Gamma\!_{f}))$) and
      $\mathbf{Dom}(\mathrm{Sg}_{\mathbf{A}\times\mathbf{D}}(\Gamma\!_{f}))$ is generated by $C$.
\item If $\mathbf{E}$ is a $C$-generated relative subalgebra of $\mathbf{A}$ and
      $g\colon \mathbf{E}\mor \mathbf{D}$ is a homomorphic extension of $f$, then
      $\Gamma\!_{g}\subseteq \Gamma\!_{f^{\flat}}$.
\item If $\mathbf{A}$ is a $\Sigma$-algebra, then
      $f^{\flat}$ is a closed homomorphism from
      $\mathbf{Dom}(\mathrm{Sg}_{\mathbf{A}\times\mathbf{D}}(\Gamma\!_{f}))$ to $\mathbf{D}$;
      and if $f^{\flat}$ is closed and $g\colon \mathbf{E}\mor \mathbf{D}$ is any closed
      homomorphic extension of $f$ such that $\mathbf{E}$ is a $C$-generated relative
      subalgebra of $\mathbf{A}$, then $g = f^{\flat}$.
\end{enumerate}
\end{lemma}

\begin{definition}
Let $\mathbf{A}$ and $\mathbf{B}$ be partial $\Sigma$-algebras and $f$ a homomorphism from $\mathbf{A}$ to $\mathbf{B}$. We will say that $f$ is \emph{dense} if $f[A]$ generates $\mathbf{B}$, i.e., if $\mathrm{Sg}_{\mathbf{B}}(f[A]) = B$.
\end{definition}

\begin{lemma}
Let $\mathbf{A}$ and $\mathbf{B}$ be partial $\Sigma$-algebras and $f$ a homomorphism from $\mathbf{A}$ to $\mathbf{B}$. Then $f$ is an epimorphism if and only if $f$ is dense.
\end{lemma}

\begin{remark}
If a homomorphism between partial $\Sigma$-algebras is dense and closed, then it is surjective (because the image of a closed homomorphism is a closed subset of the target).
\end{remark}

\begin{proposition}
Let $\mathbf{A}$ and $\mathbf{B}$ be partial $\Sigma$-algebras and $f$ a homomorphism from $\mathbf{A}$ to $\mathbf{B}$. Then there exists a unique epimorphism $f^{\mathrm{e}}$ from $\mathbf{A}$ to $\mathbf{Sg}_{\mathbf{B}}(f[A])$ such that $f = \mathrm{in}^{\mathbf{Sg}_{\mathbf{B}}[f[A]]}\circ f^{\mathrm{e}}$, where $\mathrm{in}^{\mathbf{Sg}_{\mathbf{B}}(f[A])}$ is the closed embedding of $\mathbf{Sg}_{\mathbf{B}}(f[A])$ into $\mathbf{B}$. Moreover, (1) $f^{\mathrm{e}}$ is an isomorphism if and only if $f$ is closed and injective; and (2) for every partial $\Sigma$-algebra $\mathbf{C}$, every epimorphism $g$ from $\mathbf{A}$ to $\mathbf{C}$ and every closed injective homomorphism $h$ from $\mathbf{C}$ to $\mathbf{B}$, if $f = h\circ g$, then there exists a unique isomorphism $k$ from $\mathbf{Sg}_{\mathbf{B}}(f[A])$ to $\mathbf{C}$ such that $g = k\circ f^{\mathrm{e}}$ and $h \circ k = \mathrm{in}^{\mathbf{Sg}_{\mathbf{B}}(f[A])}$.
\end{proposition}

Therefore $(\mathrm{Epimorphisms}, \mathrm{Closed\, and\, injective\, homomorphisms})$ is a factorization system in $\mathsf{PAlg}(\Sigma)$.

\section{Congruences}

Congruences provide an internal description of the homomorphic images from a many-sorted partial $\Sigma$-algebra and closed congruences provide an internal description of the closed homomorphic images from a many-sorted partial $\Sigma$-algebra. Moreover, closed congruences are fundamental for the description of the varieties of many-sorted partial $\Sigma$-algebras defined by existentially conditional existence equations.

\begin{definition}\label{DPAlgCong}
Let $\mathbf{A}$ be a partial $\Sigma$-algebra and $\Phi$ an $S$-sorted equivalence on $A$. We will say that $\Phi$ is an $S$-\emph{sorted congruence on} (or, to abbreviate, a \emph{congruence on}) $\mathbf{A}$ if, for every $(\mathbf{s},s)\in (S^{\star}-\{\lambda\})\times S$, every $\sigma\in \Sigma_{\mathbf{s},s}$,
and every $(a_{j})_{j\in\bb{\mathbf{s}}},(b_{j})_{j\in\bb{\mathbf{s}}}\in A_{\mathbf{s}}$, if, for every $j\in \bb{\mathbf{s}}$, $(a_{j}, b_{j})\in\Phi_{s_{j}}$, $(a_{j})_{j\in \bb{\mathbf{s}}}\in\mathrm{Dom}(\sigma^{\mathbf{A}})$, and $(b_{j})_{j\in \bb{\mathbf{s}}}\in\mathrm{Dom}(\sigma^{\mathbf{A}})$, then $(\sigma^{\mathbf{A}}((a_{j})_{j\in\bb{\mathbf{s}}}), \sigma^{\mathbf{A}}((b_{j})_{j\in\bb{\mathbf{s}}}))\in \Phi_{s}$.
Sometimes, we will abbreviate the expression ``$S$-sorted congruence on $\mathbf{A}$'' to ``congruence on $\mathbf{A}$''. We will denote by $\mathrm{Cgr}(\mathbf{A})$ the set of all $S$-sorted congruences on $\mathbf{A}$, which is an algebraic closure system on $A\times A$, by $\mathrm{Cg}_{\mathbf{A}}$ the corresponding algebraic closure operator, by $\mathbf{Cgr}(\mathbf{A})$ the algebraic lattice $(\mathrm{Cgr}(\mathbf{A}),\subseteq)$, by $\nabla_{\mathbf{A}}$ the greatest element of $\mathbf{Cgr}(\mathbf{A})$, and by $\Delta_{\mathbf{A}}$ the least element of $\mathbf{Cgr}(\mathbf{A})$.

For a congruence $\Phi$ on $\mathbf{A}$, the \emph{quotient partial $\Sigma$-algebra of} $\mathbf{A}$ \emph{by} $\Phi$, denoted by $\mathbf{A}/\Phi$, is the partial $\Sigma$-algebra $(A/\Phi, F^{\mathbf{A}/\Phi})$, where, for every  $(\mathbf{s},s)\in S^{\star}\times S$ and every $\sigma\in \Sigma_{\mathbf{s},s}$, the domain of the partial operation $F_{\sigma}^{\mathbf{A}/\Phi}$, also denoted, to simplify, by $\sigma^{\mathbf{A}/{\Phi}}$, from $(A/\Phi)_{\mathbf{s}}$ to $A_{s}/\Phi_{s}$ is the set
$$
\{([a_{j}]_{\Phi_{s_{j}}})_{j\in\bb{\mathbf{s}}}\in (A/\Phi)_{\mathbf{s}}\mid \exists\, (a'_{j})_{j\in\bb{\mathbf{s}}}\in A_{\mathbf{s}}\,\text{such that}\,\forall\, j\in\bb{\mathbf{s}}\,((a_{j},a'_{j})\in\Phi_{s_{j}})\},
$$
and if $([a_{j}]_{\Phi_{s_{j}}})_{j\in\bb{\mathbf{s}}}$ in $\mathrm{Dom}(\sigma^{\mathbf{A}/{\Phi}})$, then
$$
\sigma^{\mathbf{A}/{\Phi}}(([a_{j}]_{\Phi_{s_{j}}})_{j\in\bb{\mathbf{s}}}) = [\sigma^{\mathbf{A}}((a_{j})_{j\in \bb{\mathbf{s}}})]_{\Phi_{s}},
$$  %
and the \emph{canonical projection} from $\mathbf{A}$ to $\mathbf{A}/\Phi$, denoted by $\mathrm{pr}^{\Phi}\colon \mathbf{A}\mor \mathbf{A}/\Phi$, is the full and surjective homomorphism determined by the projection from $A$ to $A/\Phi$.
The ordered pair $(\mathbf{A}/\Phi,\mathrm{pr}^{\Phi})$ has the following universal property: $\mathrm{Ker}(\mathrm{pr}^{\Phi})$ is $\Phi$ and, for every $\Sigma$-algebra $\mathbf{B}$ and every homomorphism $f$ from $\mathbf{A}$ to $\mathbf{B}$, if $ \Phi\subseteq\mathrm{Ker}(f)$, then there exists a unique homomorphism $h$ from $\mathbf{A}/\Phi$ to $\mathbf{B}$ such that $h\circ\mathrm{pr}^{\Phi} = f$. In particular, if $\Psi$ is a congruence on $A$ such that $\Phi\subseteq \Psi$, then we will denote by $\mathrm{p}^{\Phi,\Psi}$ the unique homomorphism from $\mathbf{A}/\Phi$ to $\mathbf{A}/\Psi$ such that $\mathrm{p}^{\Phi,\Psi}\circ \mathrm{pr}^{\Phi} = \mathrm{pr}^\Psi$.

We will say that $\Phi$ is an $S$-\emph{sorted closed congruence on} (or, to abbreviate, a \emph{closed congruence on}) $\mathbf{A}$ if $\Phi$ is a congruence on $\mathbf{A}$ and, for every $(\mathbf{s},s)\in (S^{\star}-\{\lambda\})\times S$, every $\sigma\in \Sigma_{\mathbf{s},s}$,
and every $(a_{j})_{j\in\bb{\mathbf{s}}},(b_{j})_{j\in\bb{\mathbf{s}}}\in A_{\mathbf{s}}$, if, for every $j\in \bb{\mathbf{s}}$, $(a_{j}, b_{j})\in\Phi_{s_{j}}$, and $(a_{j})_{j\in \bb{\mathbf{s}}}\in\mathrm{Dom}(\sigma^{\mathbf{A}})$, then $(b_{j})_{j\in \bb{\mathbf{s}}}\in\mathrm{Dom}(\sigma^{\mathbf{A}})$. We will denote by $\mathrm{Cgr}_{\mathrm{c}}(\mathbf{A})$ the set of all $S$-sorted closed congruences on $\mathbf{A}$. The congruence $\Delta_{\mathbf{A}}$ is a closed congruence on $\mathbf{A}$, while $\nabla_{\mathbf{A}}$ is closed if and only if, for every $(\mathbf{s},s)\in (S^{\star}\times S$ and every $\sigma\in \Sigma_{\mathbf{s},s}$, $\sigma^{\mathbf{A}}$ is either total or discrete. If $\Phi$ and $\Psi$ are closed congruences on $\mathbf{A}$, then $\mathrm{Cg}_{\mathbf{A}}(\Phi\cup\Psi)$, their supremum in $\mathbf{Cgr}(\mathbf{A})$, is a closed congruence on $\mathbf{A}$ and it is the least equivalence relation containing them, i.e., 
$
\mathrm{Cg}_{\mathbf{A}}(\Phi\cup\Psi) = \bigcup_{n\in \mathbb{N}}(\Psi\circ\Phi)^{n}
$ 
(let us point out that, in general, the supremum of two non-closed congruences is not necessarily the least equivalence relation containing them). The set $\mathrm{Cgr}_{\mathrm{c}}(\mathbf{A})$ is a principal ideal of $\mathbf{Cgr}(\mathbf{A})$ precisely the one determined by $\bigcup_{\Phi\in \mathrm{Cgr}_{\mathrm{c}}(\mathbf{A})}\Phi$, which is the largest closed congruence on $\mathbf{A}$.
\end{definition}

\begin{remark}
Let $\mathbf{A}$ be a partial $\Sigma$-algebra and $\Phi$ an $S$-sorted equivalence on $A$. Then $\Phi$ is a congruence on $\mathbf{A}$ if and only if $\Phi$ is a subalgebra of $\mathbf{A}\times \mathbf{A}$.
\end{remark}

\begin{remark}\label{RPAlgHomFinal}
Let $\mathbf{A}$ be a partial $\Sigma$-algebra and $\Phi$ an $S$-sorted congruence on $A$, then the full and surjective homomorphism $\mathrm{pr}^{\Phi}$ from $\mathbf{A}$ to $\mathbf{A}/\Phi$ is final and surjective, where a homomorphism $f$ from a partial $\Sigma$-algebra $\mathbf{A}$ to another $\mathbf{B}$ is said to be \emph{final} if, for every partial $\Sigma$-algebra $\mathbf{C}$ and every S-sorted mapping $h$ from $B$ to $C$, if $h\circ f$ is a homomorphism from $\mathbf{A}$ to $\mathbf{C}$, then $h$ is a homomorphism from $\mathbf{B}$ to $\mathbf{C}$. For a homomorphism $f$ from a partial $\Sigma$-algebra $\mathbf{A}$ to another $\mathbf{B}$ the property of being final is very interesting because it allows one to state that an $S$-sorted mapping $h$ from $B$ to $C$, the underlying set of another partial $\Sigma$-algebra $\mathbf{C}$, is a homomorphism if its composition with the subjacent application of the homomorphism $f$ is a homomorphism from $\mathbf{A}$ to $\mathbf{C}$.
Let us point out that a homomorphism $f$ from a partial $\Sigma$-algebra $\mathbf{A} $ to another $\mathbf{B} $ is final if and only if, for every $(\mathbf{s},s)\in S^{\star}\times S$ and every $\sigma\in \Sigma_{\mathbf{s},s}$, $\Gamma_{\!\sigma^{\mathbf{B}}} = (f_{\mathbf{s}}\times f_{s})[\Gamma_{\!\sigma^{\mathbf{A}}}]$.
\end{remark}

\begin{remark}
Let $\mathbf{A}$ be a partial $\Sigma$-algebra. If $\Phi$ is an $S$-sorted closed congruence on $A$, then $\mathrm{pr}^{\Phi}\colon \mathbf{A}\mor \mathbf{A}/\Phi$, is a closed and surjective homomorphism.
\end{remark}

\begin{remark}
Let $\mathsf{ClfdPAlg}(\Sigma)$ be the category whose objects are the \emph{classified partial $\Sigma$-algebras}, i.e., the ordered pairs $(\mathbf{A},\Phi)$ where $\mathbf{A}$ is a partial $\Sigma$-algebra and $\Phi$ a congruence on $\mathbf{A}$, and in which the set of morphisms from $(\mathbf{A},\Phi)$ to $(\mathbf{B},\Psi)$ is the set of all homomorphisms $f$ from $\mathbf{A}$ to $\mathbf{B}$ such that, for every $s\in S$ and every $(x,y)\in A^{2}_{s}$, if $(x,y)\in \Phi_{s}$, then $(f_{s}(x),f_{s}(y))\in \Psi_{s}$. Let $G$ be the functor from $\mathsf{PAlg}(\Sigma)$ to $\mathsf{ClfdPAlg}(\Sigma)$ whose object mapping sends $\mathbf{A}$ to $(\mathbf{A},\Delta_{\mathbf{A}})$ and whose morphism mapping sends $f\colon \mathbf{A}\mor \mathbf{B}$ to $f\colon (\mathbf{A},\Delta_{\mathbf{A}})\mor (\mathbf{B},\Delta_{\mathbf{B}})$. Then, for every classified partial $\Sigma$-algebra $(\mathbf{A},\Phi)$, there exists a universal mapping from $(\mathbf{A},\Phi)$ to $G$, which is precisely the ordered pair $(\mathbf{A}/\Phi,\mathrm{pr}^{\Phi})$ with $\mathrm{pr}^{\Phi}\colon (\mathbf{A},\Phi)\mor (\mathbf{A}/\Phi,\Delta_{\mathbf{A}/\Phi})$.
\end{remark}

\begin{remark}\label{RPAlgLiftCo}
Let $\mathbf{A} = (A,F)$ be a partial $\Sigma$-algebra and $\Phi$ an $S$-sorted equivalence on $A$. Then there exists a unique structure of partial $\Sigma$-algebra $\mathrm{L}^{\mathrm{pr}^{\Phi}}(F)$ on $A/\Phi$, the \emph{co-optimal lift of} $F$ \emph{with respect to} $\mathrm{pr}^{\Phi}$, such that, for every partial $\Sigma$-algebra $\mathbf{C}$ and every mapping $h$ from $A/\Phi$ to $C$, if $h\circ\mathrm{pr}^{\Phi}$ is a homomorphism from $\mathbf{A}$ to, $\mathbf{C}$ then $h$ is a homomorphism from $\mathbf{C}$ to $(A/\Phi,\mathrm{L}^{\mathrm{pr}^{\Phi}}(F))$ if and only if $\Phi$ is an $S$-sorted congruence on $\mathbf{A}$. This shows that, for every $S$-sorted congruence $\Phi$ on $\mathbf{A}$ $A$, $A/\Phi$ is the underlying $S$-sorted set of exactly one $\Sigma$-algebra, and that the structure of partial $\Sigma$-algebra on $A/\Phi$ is uniquely determined by the structure of partial $\Sigma$-algebra on $A$ and the specification of the $S$-sorted congruence $\Phi$ on $A$. What we have just said is, in fact, a particular case of the following: for a partial $\Sigma$-algebra $\mathbf{A} = (A,F)$ and an $S$-sorted mapping $f$ from $A$ to $B$, there exists a co-optimal lift of $F$ with respect to $f$ if and only if $\mathrm{Ker}(f)$ is a congruence on $\mathbf{A}$.
\end{remark}

We next state the weak homomorphis theorem.

\begin{proposition}
Let $\mathbf{A}$, $\mathbf{B}$, and $\mathbf{C}$ be partial $\Sigma$-algebras, $f$ a full and surjective homomorphism from $\mathbf{A}$ to $\mathbf{B}$, and $g$ a homomorphism from $\mathbf{A}$ to $\mathbf{C}$.
If $\mathrm{Ker}(g)\supseteq \mathrm{Ker}(f)$, then there exists a unique homomorphism $h$ from $\mathbf{B}$ to $\mathbf{C}$ such that $g = h\circ f$ (obviously, if there exists a homomorphism $h$ from $\mathbf{B}$ to $\mathbf{C}$ such that $g = h\circ f$, then  $\mathrm{Ker}(g)\supseteq \mathrm{Ker}(f)$). Moreover, (1) $h$ is a full homomorphism if and only if $g$ is a full homomorphism; (2) $h$ is an isomorphism if and only if $g$ is a surjective and full homomorphism and $\mathrm{Ker}(g) = \mathrm{Ker}(f)$; (3) if $g$ is closed, then $h$ is closed; (4) if $f$ is closed, then $h$ is closed if and only if $g$ is closed.
\end{proposition}

From the just stated proposition we obtain the weak isomorphism theorem.

\begin{corollary}
Let $\mathbf{A}$, $\mathbf{B}$ be partial $\Sigma$-algebras, and $f$ a full and surjective homomorphism from $\mathbf{A}$ to $\mathbf{B}$. Then $\mathbf{B}$ and $\mathbf{A}/\mathrm{Ker}(f)$ are isomorphic.
\end{corollary}

\begin{proposition}
Let $\mathbf{A}$ and $\mathbf{B}$ be partial $\Sigma$-algebras and $f$ a homomorphism from $\mathbf{A}$ to $\mathbf{B}$. Then there exists a unique injective homomorphism $f^{\mathrm{m}}$ from $\mathbf{A}/\mathrm{Ker}(f)$ to $\mathbf{B}$ such that $f = f^{\mathrm{m}}\circ \mathrm{pr}^{\mathrm{Ker}(f)}$. Moreover, (1) $f^{\mathrm{m}}$ is an isomorphism if and only if $f$ is full and surjective; and (2) for every partial $\Sigma$-algebra $\mathbf{C}$, every full and surjective homomorphism $g$ from $\mathbf{A}$ to $\mathbf{C}$ and every injective homomorphism $h$ from $\mathbf{C}$ to $\mathbf{B}$, if $f = h\circ g$, then there exists a unique isomorphism $k$ from $\mathbf{C}$ to $\mathbf{A}/\mathrm{Ker}(f)$ such that $k\circ g =  \mathrm{pr}^{\mathrm{Ker}(f)}$ and $f^{\mathrm{m}}\circ k = h$.
\end{proposition}

Therefore $(\mathrm{Full\, and\, surjective\, homomorphisms}, \mathrm{Monomorphisms})$ is a factorization system in $\mathsf{PAlg}(\Sigma)$.

\section{Partial Dedekind-Peano Algebras}

We next define for many-sorted partial algebras the counterpart of the notion of Dedekind-Peano algebra. Moreover, we state the principle of definition by algebraic recursion for free many-sorted $\Sigma$-algebras, i.e., many-sorted DP-algebras, with respect to many-sorted partial $\Sigma$-algebras, which, for many-sorted partial algebras, is more important than the principle of definition by algebraic recursion for partial  Dedekind-Peano algebra with respect to many-sorted algebras.

\begin{definition}\label{DPAlgDP}
Let $\mathbf{A}$ be a partial $\Sigma$-algebra. We will say that $\mathbf{A}$ is a \emph{partial Dedekind-Peano $\Sigma$-algebra}, abbreviated to $\mathrm{PDP}$-\emph{algebra} when this is unlikely to cause confusion, if the following axioms hold
\begin{enumerate}
\item[PDP1.] For every $(\mathbf{s},s)\in S^{\star}\times S$ and every $\sigma\in\Sigma_{\mathbf{s},s}$, $\sigma^{\mathbf{A}}\colon A_{\mathbf{s}}\dmor A_{s}$ is injective, i.e., for every $(a_{j})_{j\in\bb{\mathbf{s}}}$, $(b_{j})_{j\in\bb{\mathbf{s}}}\in \mathrm{Dom}(\sigma^{\mathbf{A}})$, if $\sigma^{\mathbf{A}}((a_{j})_{j\in\bb{\mathbf{s}}}) = \sigma^{\mathbf{A}}((b_{j})_{j\in\bb{\mathbf{s}}})$, then $(a_{j})_{j\in\bb{\mathbf{s}}}=(b_{j})_{j\in\bb{\mathbf{s}}}$.

\item[PDP2.] For every $s\in S $ and every $\sigma, \tau \in\Sigma_{\cdot,s}$, if $\sigma\neq\tau$, then $\mathrm{Im}(\sigma^{\mathbf{A}})\cap \mathrm{Im}(\tau^{\mathbf{A}}) = \varnothing$, i.e., $\sigma^{\mathbf{A}}[\mathrm{Dom}(\sigma^{\mathbf{A}})]\cap \tau^{\mathbf{A}}[\mathrm{Dom}(\tau^{\mathbf{A}})] = \varnothing$.

\item[PDP3.] $\mathrm{Sg}_{\mathbf{A}}(A-(\bigcup_{\sigma\in\Sigma_{\cdot,s}}\mathrm{Im}(\sigma^{\mathbf{A}}))_{s\in S}) = A$.
\end{enumerate}
We call the $S$-sorted set $A-(\bigcup_{\sigma\in\Sigma_{\cdot,s}}\mathrm{Im}(\sigma^{\mathbf{A}}))_{s\in S}$
the \emph{basis of Dedekind-Peano of} $\mathbf{A}$, and we denote it by $\mathrm{B}(\mathbf{A})$.
\end{definition}

\begin{remark}
Every partial Dedekind-Peano algebra can be obtained, up to isomorphism, as a weak subalgebra of a Dedekind-Peano algebra. But unlike for many-sorted algebras in a variety, for many-sorted partial algebras what we have is that every many-sorted partial algebra in an $\mathrm{E}$-variety is a \emph{weak} homomorphic image of a free many-sorted partial algebra, yet not necessarily a \emph{full} or \emph{closed} homomorphic image. And in $\mathrm{ECE}$-varieties or $\mathrm{QE}$-varieties the more general concept of universal solution becomes of more importance than the one of freeness. The notions of $\mathrm{E}$-, $\mathrm{ECE}$-, and $\mathrm{QE}$-variety will be defined later on.
\end{remark}

\begin{remark}
A discrete partial $\Sigma$-algebra is always a PDP-algebra on its underlying $S$-sorted set. Moreover, every weak subalgebra of a PDP-algebra is a PDP-algebra.
\end{remark}

We next state the principle of the definition by algebraic recursion for many-sorted PDP-algebras with respect to many-sorted algebras.

\begin{proposition}\label{PPAlgDP}
Let $\mathbf{A}$ be a PDP-algebra, $\mathrm{B}(\mathbf{A})$ its basis, $\mathbf{B}$ a $\Sigma$-algebra, and $f$ a mapping from $\mathrm{B}(\mathbf{A})$ to $B$. Then there exists a unique homomorphism $f^{\sharp}$ from $\mathbf{A}$ to $\mathbf{B}$ such that $f^{\sharp}\circ \mathrm{in}^{\mathrm{B}(\mathbf{A}),\mathbf{A}} = f$.
\end{proposition}

\begin{proof}
It suffices to take as $f^{\sharp}$ precisely $(\mathbf{A},\Gamma\!_{f^{\sharp}},\mathbf{B})$, where $\Gamma\!_{f^{\sharp}}$, the underlying $S$-sorted function of $f^{\sharp}$ is $\mathrm{Sg}_{\mathbf{A}\times \mathbf{B}}(\Gamma\!_{f})$.
\end{proof}

We next state the principle of definition by algebraic recursion for free many-sorted $\Sigma$-algebras, i.e., many-sorted DP-algebras, with respect to many-sorted partial $\Sigma$-algebras.

\begin{proposition}[Proposition~30, pp.~29--32,~\cite{CVCL23}]\label{PPAlgPartial}
Let $X$ be an $S$-sorted set, $\mathbf{A}$ a partial $\Sigma$-algebra and $f$ an $S$-sorted mapping from $X$ to $A$. Then there exists a unique homomorphism 
$$
f^{\partial}\colon \boldsymbol{\partial}(f)\mor \mathbf{A}
$$ 
such that
\begin{enumerate}
\item $\eta^{X}[X]\subseteq \partial(f)$, where $\eta^{X}$ is the canonical embedding of $X$ into $\mathrm{T}_{\Sigma}(X)$.
\item $\boldsymbol{\partial}(f)$ is an $X(\cong\eta^{X}[X])$-generated relative subalgebra of $\mathbf{T}_{\Sigma}(X)$.
\item $f^{\partial}\circ \mathrm{in}^{X,\boldsymbol{\partial}(f)} = f$, where $\mathrm{in}^{X,\boldsymbol{\partial}(f)}$ is the canonical embedding of $X$ into $\boldsymbol{\partial}(f)$.
\item $f^{\partial}$ is a closed homomorphism.
\item $f^{\partial}$ is the largest homomorphic extension of $f$ to an $X$-generated relative subalgebra of $\mathbf{T}_{\Sigma}(X)$ with codomain $\mathbf{A}$.
\end{enumerate}
\end{proposition}
\begin{proof}
Let $\mathbf{A}^{\infty}$ be the $\Sigma$-algebra defined as follows: The underlying $S$-sorted set of $\mathbf{A}^{\infty}$ is $(A_{s}\cup\{A_{s}\})_{s\in S}$ and, for every $(\mathbf{s},s)\in S^{\star}\times S$ and every $\sigma\in \Sigma_{\mathbf{s},s}$, the operation $\sigma^{\mathbf{A}^{\infty}}\colon A^{\infty}_{\mathbf{s}}\mor
A^{\infty}_{s}$ associated to $\sigma$ is defined as
$$
\sigma^{\mathbf{A}^{\infty}}
\nfunction
{A^{\infty}_{\mathbf{s}}}{A^{\infty}_{s}}
{(a_{j})_{j\in \bb{\mathbf{s}}}}{
\begin{cases}
\sigma^{\mathbf{A}}((a_{j})_{j\in \bb{\mathbf{s}}}), & \text{if } (a_{j})_{j\in \bb{\mathbf{s}}}\in \mathrm{Dom}(\sigma^{\mathbf{A}});\\
  A_{s}, & \text{otherwise.}
\end{cases}
   }
$$

We will call $\mathbf{A}^{\infty}$ the \emph{one-point per sort completion} of the partial $\Sigma$-algebra $\mathbf{A}$. Notice that $\mathbf{A}$ is a relative subalgebra of $\mathbf{A}^{\infty}$, that $\mathbf{A}^{\infty}$ need not be generated by $A$ (since some of the new elements are not necessarily accessible from $A$) and that the canonical embedding $\mathrm{in}^{A,A^{\infty}}$ of $A$ into $A^{\infty}$ is the underlying mapping of a full and injective homomorphism, denoted by $\mathrm{in}^{\mathbf{A},\mathbf{A}^{\infty}}$, from $\mathbf{A}$ to $\mathbf{A}^{\infty}$.
Then, by the universal property of the free algebra, there exists a unique homomorphism 
$$
(\mathrm{in}^{\mathbf{A},\mathbf{A}^{\infty}}\circ f)^{\sharp}\colon\mathbf{T}_{\Sigma}(X)\mor\mathbf{A}^{\infty}
$$ 
such that $(\mathrm{in}^{\mathbf{A},\mathbf{A}^{\infty}}\circ f)^{\sharp}\circ\eta^{X} = \mathrm{in}^{\mathbf{A},\mathbf{A}^{\infty}}\circ f$. Let $\partial(f)$ be $((\mathrm{in}^{A,A^{\infty}}\circ f)^{\sharp})^{-1}[A]$ and let $f^{\partial}$ be $(\mathrm{in}^{A,A^{\infty}}\circ f)^{\sharp}\bigr|_{((\mathrm{in}^{A,A^{\infty}}\circ f)^{\sharp})^{-1}[A]}^{A}$, i.e., the birestriction of $(\mathrm{in}^{A,A^{\infty}}\circ f)^{\sharp}$ to $((\mathrm{in}^{A,A^{\infty}}\circ f)^{\sharp})^{-1}[A]$ and $A$. Then $\boldsymbol{\partial}(f)$, the relative subalgebra of $\mathbf{T}_{\Sigma}(X)$ on $\partial(f)$, together with the homomorphism determined by  $f^{\partial}$, which, with the customary abuse of notation, we denote by the same symbol, satisfy the desired conditions.

To better understand the mappings at play, we provide the commutative diagram in Figure~\ref{FPAlgPartial}.
\end{proof}

\begin{figure}
$$\xymatrix@C=60pt@R=40pt{
    X
    \ar@/_1pc/[ddr]_{f}
    \ar@/^1pc/[drr]^{\eta^{X}}
    \ar[dr]|-{\mathrm{in}^{X,\boldsymbol{\partial}(f)}} \\
     & \boldsymbol{\partial}(f)
     \ar[d]_{f^{\partial}}
     \ar[r]^-{\mathrm{in}^{\boldsymbol{\partial}(f),\mathbf{T}_{\Sigma}(X)}}  &
     \mathbf{T}_{\Sigma}(X)
     \ar[d]^{(\mathrm{in}^{\mathbf{A},\mathbf{A}^{\infty}}\circ f)^{\sharp}} \\
     &
     \mathbf{A}
     \ar[r]_{\mathrm{in}^{\mathbf{A},\mathbf{A}^{\infty}}} &
     \mathbf{A}^{\infty}
   }
$$
\caption{Recursion Theorem for DP-algebras with respect to partial algebras.}
\label{FPAlgPartial}
\end{figure}

\begin{remark}
Let us note that, by the Axiom of Regularity, we have that, for every $s\in S$, $A_{s}\cap\{A_{s}\} = \varnothing$. Moreover, for every $(\mathbf{s},s)\in S^{\star}\times S$ and every $\sigma\in \Sigma_{\mathbf{s},s}$, $\Gamma_{\sigma^{\boldsymbol{\partial}(f)}}$, the underlying partial function of $\sigma^{\boldsymbol{\partial}(f)}$, the partial operation from $\partial(f)_{\mathbf{s}}$ to $\partial(f)_{s}$ associated to $\sigma$, is
$$
\Gamma_{\!\sigma^{\mathbf{T}_{\Sigma}(X)}}\cap (\partial(f)_{\mathbf{s}}\times \partial(f)_{s}).
$$
On the other hand, $\Gamma\!_{f^{\partial}} = \mathrm{Sg}_{\mathbf{T}_{\Sigma}(X)\times \mathbf{A}}((\eta^{X}\times \mathrm{id}^{A})[\Gamma\!_{f}])$, i.e., the underlying function of $f^{\partial}$, is the subalgebra of the partial $\Sigma$-algebra $\mathbf{T}_{\Sigma}(X)\times \mathbf{A}$ generated by the image under the $S$-sorted mapping $\eta^{X}\times \mathrm{id}^{A}\colon X\times A\mor\mathrm{T}_{\Sigma}(X)\times \mathrm{A}$, of the underlying function of $f$. Finally, we have that
\begin{enumerate}
\item For every $s\in S$ and every $x\in X_{s}$, $f^{\partial}_{s}(x) = f_{s}(x)$.
\item For every $s\in S$ and every $\sigma\in\Sigma_{\lambda,s}$, if $\mathrm{Dom}(\sigma^{\mathbf{A}})\neq\varnothing$, then $\sigma\in \partial(f)_{s}$ and $f^{\partial}_{s}(\sigma^{\boldsymbol{\partial}(f)}) = \sigma^{\mathbf{A}}$.
\item For every $(\mathbf{s},s)\in (S^{\star}-\{\lambda\})\times S$, every $\sigma\in \Sigma_{\mathbf{s},s}$, and every $(P_{j})_{j\in \bb{\mathbf{s}}}\in \mathrm{T}_{\Sigma}(X)_{\mathbf{s}}$, if, for every $j\in\bb{\mathbf{s}}$, $P_{j}\in \partial(f)_{s_{j}}$ and the fami\-ly $(f^{\partial}_{s_{j}}(P_{j}))_{j\in \bb{\mathbf{s}}}\in \mathrm{Dom}(\sigma^{\mathbf{A}})$, then $\sigma^{\boldsymbol{\partial}(f)}((P_{j})_{j\in \bb{\mathbf{s}}})\in \partial(f)_{s}$ and it holds that $$f^{\partial}_{s}(\sigma^{\boldsymbol{\partial}(f)}((P_{j})_{j\in \bb{\mathbf{s}}})) = \sigma^{\mathbf{A}}((f^{\partial}_{s_{j}}(P_{j})))_{j\in \bb{\mathbf{s}}}).$$
\end{enumerate}
\end{remark}

\begin{remark}
The construction that assigns to a partial $\Sigma$-algebra $\mathbf{A}$ the $\Sigma$-algebra $\mathbf{A}^{\infty}$ is not the object mapping of a functor left adjoint to the inclusion functor $\mathrm{In}_{\mathsf{Alg}(\Sigma)}$ from $\mathsf{Alg}(\Sigma)$ to $\mathsf{PAlg}(\Sigma)$. However, the just mentioned construction is the object mapping of a functor $(\cdot)^{\infty}$ from the category $\mathsf{PAlg}(\Sigma)_{\mathrm{c}}$ of many-sorted partial $\Sigma$-algebras and closed homomorphisms to the category $\mathsf{Alg}(\Sigma)$.
\end{remark}

\section{
\texorpdfstring
{The free completion of a many-sorted partial $\Sigma$-algebra}
{The free completion of a many-sorted partial algebra}
}

We next prove, as a consequence of a well-known theorem about adjoint functors, that the inclusion functor $\mathrm{In}_{\mathsf{Alg}(\Sigma)}$ from $\mathsf{Alg}(\Sigma)$ to $\mathsf{PAlg}(\Sigma)$ has a left adjoint, the (absolutely) free completion functor. We point out that the free completion of a many-sorted partial $\Sigma$-algebra is one of the most useful tools of the theory of partial algebras. At the end of this subsection, and related to the free completion functor, we state the generalized recursion theorem of Schmidt (see~\cite{sch70}), which we will also call the Schmidt construction. This construction is fundamental in the field of many-sorted partial algebras.

\begin{figure}
$$\xymatrix{
   \mathsf{A}
   \ar[r]^{G}
   \ar[dr]_{U}
   & \mathsf{B} \ar[d]^{V} \\
   & \mathsf{C}
    }
$$
\caption{Conditions for the existence of a left adjoint.}
\label{FPAlgCondAdj}
\end{figure}

\begin{proposition}
If the diagram of categories and functors in Figure~\ref{FPAlgCondAdj}
commutes, and the following conditions are satisfied:
\begin{enumerate}
\item[(i)] $\mathsf{A}$ is complete, well-powered and co-well-powered,
\item[(ii)] $G$ preserves limits,
\item[(iii)] $U$ has a left-adjoint,
\item[(iv)] $V$ is faithful,
\end{enumerate}
then $G$ has a left adjoint.
\end{proposition}

\begin{figure}
$$\xymatrix{
   \mathsf{Alg}(\Sigma)
   \ar[r]^{\mathrm{In}_{\mathsf{Alg}(\Sigma)}}
   \ar[dr]_{G_{\mathsf{Alg}(\Sigma)}}
   & \mathsf{PAlg}(\Sigma) \ar[d]^{G_{\mathsf{PAlg}(\Sigma)}} \\
   & \mathsf{Set}^{S}
    }
$$
\caption{Conditions for the existence of the free completion.}
\label{FFreeComp}
\end{figure}

\begin{corollary}
The diagram in Figure~\ref{FFreeComp}
commutes and
\begin{enumerate}
\item[(i)] $\mathsf{Alg}(\Sigma)$ is complete, well-powered and co-well-powered,
\item[(ii)] $\mathrm{In}_{\mathsf{Alg}(\Sigma)}$ preserves limits,
\item[(iii)] $G_{\mathsf{Alg}(\Sigma)}$ has a left-adjoint ($\mathbf{T}_{\Sigma}$),
\item[(iv)] $G_{\mathsf{PAlg}(\Sigma)}$ is faithful.
\end{enumerate}
Therefore $\mathrm{In}_{\mathsf{Alg}(\Sigma)}$ has a left adjoint $\mathbf{F}_{\Sigma}$, the \emph{free completion} functor.
\end{corollary}

We next provide an explicit construction of the free completion of a partial $\Sigma$-algebra.


\begin{proposition}\label{PFreeComp}
Let $\mathbf{A}$ be a partial $\Sigma$-algebra. Then there exists a $\Sigma$-algebra $\mathbf{F}_{\Sigma}(\mathbf{A})$, the free completion of $\mathbf{A}$, and a homomorphism $\eta^{\mathbf{A}}$ from $\mathbf{A}$ to $\mathbf{F}_{\Sigma}(\mathbf{A})$ such that, for every $\Sigma$-algebra $\mathbf{B}$ and every homomorphism $f$ from $\mathbf{A}$ to $\mathbf{B}$, there exists a unique homomorphism $f^{\mathrm{fc}}$ from $\mathbf{F}_{\Sigma}(\mathbf{A})$ to $\mathbf{B}$ such that $f = f^{\mathrm{fc}}\circ \eta^{\mathbf{A}}$.
\end{proposition}

\begin{proof}
Let $\mathbf{T}_{\Sigma}(A)$ be the free $\Sigma$-algebra on $A$ and, for every $(\mathbf{s},s)\in S^{\star}\times S$ and every $\sigma\in \Sigma_{\mathbf{s},s}$, let $\overline{F}_{\sigma}$ be the mapping from $\mathrm{T}_{\Sigma}(A)_{\mathbf{s}}$ to $\mathrm{T}_{\Sigma}(A)_{s}$ defined as:
$$
\overline{F}_{\sigma}
\nfunction
{\mathrm{T}_{\Sigma}(A)_{\mathbf{s}}}{\mathrm{T}_{\Sigma}(A)_{s}}
{(P_{j})_{j\in \bb{\mathbf{s}}}}{
\begin{cases}
 F^{\mathbf{A}}_{\sigma}((P_{j})_{j\in \bb{\mathbf{s}}}), & \text{if } (P_{j})_{j\in \bb{\mathbf{s}}}\in \mathrm{Dom}(F^{\mathbf{A}}_{\sigma});\\
 F^{\mathbf{T}_{\Sigma}(A)}_{\sigma}((P_{j})_{j\in \bb{\mathbf{s}}}), & \text{otherwise.}
\end{cases}
   }
$$
Let $\overline{\mathbf{T}}_{\Sigma}(\mathbf{A}) = (\mathrm{T}_{\Sigma}(A),\overline{F})$ be the resulting $\Sigma$-algebra. Then 
$\mathbf{F}_{\Sigma}(\mathbf{A}) = (\mathrm{F}_{\Sigma}(\mathbf{A}),\overline{F}) = \mathbf{Sg}_{\overline{\mathbf{T}}_{\Sigma}(\mathbf{A})}(A)$, the subalgebra of $\overline{\mathbf{T}}_{\Sigma}(\mathbf{A})$ generated by $A$, together with $\eta^{\mathbf{A}}$, the canonical embedding of $\mathbf{A}$ into $\mathbf{F}_{\Sigma}(\mathbf{A})$ induced by the corestriction of the mapping $\eta^{A}\colon A\mor \mathrm{T}_{\Sigma}(A)$ to the underlying $S$-sorted set of $\mathbf{F}_{\Sigma}(\mathbf{A})$, is a universal morphism from $\mathbf{A}$ to $\mathrm{In}_{\mathbf{Alg}(\Sigma)}$. In fact, from $\mathbf{F}_{\Sigma}(\mathbf{A})$ we obtain a partial Dedekind-Peano $\Sigma$-algebra $\mathbf{F}^{\ast}_{\Sigma}(\mathbf{A}) = (\mathrm{F}_{\Sigma}(\mathbf{A}),F^{\ast})$ by defining, for every $(\mathbf{s},s)\in S^{\star}\times S$ and every $\sigma\in \Sigma_{\mathbf{s},s}$, $F_{\sigma}^{\ast}$ to be the partial mapping from $\mathrm{F}_{\Sigma}(\mathbf{A})_{\mathbf{s}}$ to $\mathrm{F}_{\Sigma}(\mathbf{A})_{s}$ whose domain of definition is $\mathrm{F}_{\Sigma}(\mathbf{A})_{\mathbf{s}}-\mathrm{Dom}(F_{\sigma}^{\mathbf{A}})$ and is such that, 
for $(P_{j})_{j\in \bb{\mathbf{s}}}\in \mathrm{Dom}(F_{\sigma}^{\ast})$, 
$F_{\sigma}^{\ast}((P_{j})_{j\in \bb{\mathbf{s}}}) = \overline{F}_{\sigma}((P_{j})_{j\in \bb{\mathbf{s}}})$. Then, given a $\Sigma$-algebra $\mathbf{B}$ and a homomorphism $f$ from $\mathbf{A}$ to $\mathbf{B}$, by Proposition~\ref{PPAlgDP}, 
there exists a unique homomorphism $f^{\sharp}$ from $\mathbf{F}^{\ast}_{\Sigma}(\mathbf{A})$ to $\mathbf{B}$ such that $f^{\sharp}\circ\eta^{\mathbf{A}} = f$. We recall that 
$\Gamma\!_{f^{\sharp}} = \mathrm{Sg}_{\mathbf{F}^{\ast}_{\Sigma}(\mathbf{A})\times \mathbf{B}}(\Gamma_{f})$. Since $f^{\sharp}\circ\eta^{\mathbf{A}} = f$ is a homomorphism from $\mathbf{A}$ to $\mathbf{B}$, $f^{\sharp}$ is a homomorphism from $\mathbf{F}_{\Sigma}(\mathbf{A})$ to $\mathbf{B}$. Then it suffices to take as  $f^{\mathsf{fc}}$ the homomorphism $f^{\sharp}$ but considered as a homomorphism from $\mathbf{F}_{\Sigma}(\mathbf{A})$ to $\mathbf{B}$, i.e., $\Gamma\!_{f^{\mathsf{fc}}} = \Gamma\!_{f^{\sharp}}$ and $\mathrm{d}_{0}(f^{\mathsf{fc}}) = \mathbf{F}_{\Sigma}(\mathbf{A})$.
\end{proof}

\begin{corollary} \label{CFreeAdj}
The functor $\mathbf{F}_{\Sigma}$ from $\mathsf{PAlg}(\Sigma)$ to $\mathsf{Alg}(\Sigma)$, which sends a partial $\Sigma$-algebra $\mathbf{A}$ to $\mathbf{F}_{\Sigma}(\mathbf{A})$ and a homomorphism $f$ from $\mathbf{A}$ to $\mathbf{B}$ to $f^{@}$, the unique homomorphism $(\eta^{\mathbf{A}}\circ f)^{\mathrm{fc}}$ from $\mathbf{F}_{\Sigma}(\mathbf{A})$ to $\mathbf{F}_{\Sigma}(\mathbf{B})$ such that $(\eta^{\mathbf{B}}\circ f)^{\mathrm{fc}}\circ\eta^{\mathbf{A}} = \eta^{\mathbf{B}}\circ f$, is a left adjoint of the functor $\mathrm{In}_{\mathsf{Alg}(\Sigma)}$ from $\mathsf{Alg}(\Sigma)$ to $\mathsf{PAlg}(\Sigma)$.
We will call the $\Sigma$-homomorphism  $f^{@}$ the \emph{free completion of} $f$.
\end{corollary}

\begin{remark}
For every $\Sigma$-algebra $\mathbf{A}$, $(\mathbf{A},\mathrm{id}^{\mathbf{A}})$ is a free completion of 
$\mathbf{A}$, since it is a universal morphism from $\mathbf{A}$ to 
$\mathrm{In}_{\mathsf{Alg}(\Sigma)}$. 
\end{remark}

\begin{remark}
Since the functor $\mathrm{In}_{\mathsf{Alg}(\Sigma)}$ from $\mathsf{Alg}(\Sigma)$ to $\mathsf{PAlg}(\Sigma)$ is, clearly, full, faithful and injective on objects, we can assert that the functor $\mathbf{F}_{\Sigma}$ from $\mathsf{PAlg}(\Sigma)$ to $\mathsf{Alg}(\Sigma)$ is, not only a left adjoint, but also a left inverse of $\mathrm{In}_{\mathsf{Alg}(\Sigma)}$, i.e., the counit of the adjunction is the identity. We recall that $\eta$, the unit of the adjunction, determines, for every partial $\Sigma$-algebra $\mathbf{A}$ and every $\Sigma$-algebra $\mathbf{B}$ a natural isomorphism $\varphi_{\mathbf{A},\mathbf{B}}$ from 
$\mathrm{Hom}(\mathbf{F}_{\Sigma}(\mathbf{A}),\mathbf{B})$ to 
$\mathrm{Hom}(\mathbf{A},\mathrm{In}_{\mathsf{Alg}(\Sigma)}(\mathbf{B}))$ and that $\varepsilon^{\mathbf{A}}$, the value of the counit 
$\varepsilon$ of the adjunction at a $\Sigma$-algebra $\mathbf{A}$, is 
$\varphi^{-1}_{\mathrm{In}_{\mathsf{Alg}(\Sigma)}(\mathbf{A}),\mathbf{A}}(\mathrm{id}^{\mathbf{A}})$, which in this case is, precisely, $\mathrm{id}^{\mathbf{A}}$.
\end{remark}

\begin{remark}
For every partial $\Sigma$-algebra $\mathbf{A}$, the injective homomorphism $\eta^{\mathbf{A}}$ from $\mathbf{A}$ to $\mathbf{F}_{\Sigma}(\mathbf{A})$, i.e., the value of the unit of the adjuction at $\mathbf{A}$, is, in addition, an epimorphism, hence a bimorphism. Therefore, the isomorphism-closed, full subcategory $\mathsf{Alg}(\Sigma)$ of $\mathsf{PAlg}(\Sigma)$ is monoreflective and epireflective.
\end{remark}

\begin{figure}
\begin{center}
\begin{tikzpicture}
[ACliment/.style={-{To [angle'=45, length=5.75pt, width=4pt, round]}
}, scale=0.7]
\node[] (A) 	at (0,-1) 	[] {$\mathbf{A}$};
\node[] (B) 	at (0,-4) 	[] {$\mathbf{B}$};
\node[] (TA) 	at (10,2) 	[] {$\mathbf{T}_{\Sigma}(\mathbf{A})$};
\node[] (TB) 	at (10,-7) 	[] {$\mathbf{T}_{\Sigma}(\mathbf{B})$};
\node[] (TpA) 	at (7,.5) 	[] 
{$\overline{\mathbf{T}}_{\Sigma}(\mathbf{A})$};
\node[] (TpB) 	at (7,-5.5) 	[] 
{$\overline{\mathbf{T}}_{\Sigma}(\mathbf{B})$};
\node[] (FA) 	at (4,-1) 	[] {$\mathbf{F}_{\Sigma}(\mathbf{A})$};
\node[] (FB) 	at (4,-4) 	[] {$\mathbf{F}_{\Sigma}(\mathbf{B})$};
\draw[ACliment]  (A) 	to node [left]	{$f$} (B);
\draw[ACliment]  (TA) 	to node [right]	{$f^{@} = (\eta^{\mathbf{B}}\circ f)^{\sharp}$} (TB);
\draw[ACliment, bend left]  (A) 	to node [above]	{$\eta^{\mathbf{A}}$} (TA);
\draw[ACliment, bend right]  (B) 	to node [below]	{$\eta^{\mathbf{B}}$} (TB);
\draw[ACliment, bend left=20]  (A) 	to node [above]	{$\mathrm{in}^{\mathbf{A}}$} (TpA);
\draw[ACliment, bend right=20]  (B) 	to node [below]	{$\mathrm{in}^{\mathbf{B}}$} (TpB);
\draw[ACliment]  (TA) 	to node [above left]	{$(\mathrm{in}^{\mathbf{A}})^{\sharp}$} (TpA);
\draw[ACliment]  (TB) 	to node [below left]	{$(\mathrm{in}^{\mathbf{B}})^{\sharp}$} (TpB);
\draw[ACliment]  (A) 	to node [below]	{$\eta^{\mathbf{A}}$} (FA);
\draw[ACliment]  (B) 	to node [above]	{$\eta^{\mathbf{B}}$} (FB);
\draw[ACliment]  (FA) 	to node [right]	{$f^{@} = (\eta^{\mathbf{B}}\circ f)^{\mathrm{fc}}$} (FB);
\draw[ACliment]  (FA) 	to node [above left]	{$\mathrm{in}^{\mathbf{F}_{\Sigma}(\mathbf{A})}$} (TpA);
\draw[ACliment]  (FB) 	to node [below left]	{$\mathrm{in}^{\mathbf{F}_{\Sigma}(\mathbf{B})}$} (TpB);
\draw[ACliment, bend left]  (TA) 	to node [below right]	{$(\eta^{\mathbf{A}})^{\sharp}$} (FA);
\draw[ACliment, bend right]  (TB) 	to node [above right]	{$(\eta^{\mathbf{B}})^{\sharp}$} (FB);
\end{tikzpicture}
\end{center}
\caption{The free completion functor.}
\label{FFreeCompFun}
\end{figure}

\begin{remark}\label{RFreeDisc}
Let $X$ be an $S$-sorted set. Then, for $\mathbf{D}_{\Sigma}(X)$, the discrete many-sorted partial $\Sigma$-algebra associated to an $S$-sorted set $X$, the three many-sorted $\Sigma$-algebras $\mathbf{F}_{\Sigma}(\mathbf{D}_{\Sigma}(X))$, $\overline{\mathbf{T}}_{\Sigma}(\mathbf{D}_{\Sigma}(X))$, and $\mathbf{T}_{\Sigma}(\mathbf{D}_{\Sigma}(X))$ are  naturally isomorphic to $\mathbf{T}_{\Sigma}(X)$, the free many-sorted $\Sigma$-algebra on $X$. This follows from the fact that 
$\mathbf{D}_{\Sigma}\dashv G_{\mathsf{PAlg}(\Sigma)}$ (see Remark~\ref{RDiscAdj}), 
$\mathbf{F}_{\Sigma}\dashv \mathrm{In}_{\mathsf{Alg}(\Sigma)}$ and 
$G_{\mathsf{PAlg}(\Sigma)}\circ \mathrm{In}_{\mathsf{Alg}(\Sigma)} = G_{\mathsf{Alg}(\Sigma)}$.
\end{remark}

To better understand the working of the functor $\mathbf{F}_{\Sigma}$, we provide the commutative diagram in Figure~\ref{FFreeCompFun}.

Let us note that for $f\colon \mathbf{D}_{\Sigma}(X)\mor \mathbf{B}$ the diagram in Figure~\ref{FFreeCompFun} becomes the commutative diagram in Figure~\ref{FFreeCompFunDisc}

\begin{figure}
\begin{center}
\begin{tikzpicture}
[ACliment/.style={-{To [angle'=45, length=5.75pt, width=4pt, round]}
}, scale=0.7]
\node[] (A) 	at (0,-1) 	[] {$\mathbf{D}_{\Sigma}(X)$};
\node[] (B) 	at (0,-4) 	[] {$\mathbf{B}$};
\node[] (TA) 	at (4,-1) 	[] {$\mathbf{T}_{\Sigma}(X)$};
\node[] (TB) 	at (10,-7) 	[] {$\mathbf{T}_{\Sigma}(\mathbf{B})$};
\node[] (TpB) 	at (7,-5.5) 	[] 
{$\overline{\mathbf{T}}_{\Sigma}(\mathbf{B})$};
\node[] (FB) 	at (4,-4) 	[] {$\mathbf{F}_{\Sigma}(\mathbf{B})$};
\draw[ACliment]  (A) 	to node [left]	{$f$} (B);
\draw[ACliment, bend left=45]  (TA) 	to node [above right]{$f^{@}  = (\eta^{\mathbf{B}}\circ f)^{\sharp}$} (TB);
\draw[ACliment, bend right]  (B) 	to node [below]	{$\eta^{\mathbf{B}}$} (TB);
\draw[ACliment, bend right=20]  (B) 	to node [below]	{$\mathrm{in}^{\mathbf{B}}$} (TpB);
\draw[ACliment]  (TB) 	to node [below left]	{$(\mathrm{in}^{\mathbf{B}})^{\sharp}$} (TpB);
\draw[ACliment]  (A) 	to node [below]	{$\eta^{\mathbf{A}}$} (TA);
\draw[ACliment]  (B) 	to node [above]	{$\eta^{\mathbf{B}}$} (FB);
\draw[ACliment]  (TA) 	to node [right]	{$f^{@}  = (\eta^{\mathbf{B}}\circ f)^{\mathrm{fc}}$} (FB);
\draw[ACliment]  (FB) 	to node [below left]	{$\mathrm{in}^{\mathbf{F}_{\Sigma}(\mathbf{B})}$} (TpB);
\draw[ACliment, bend right]  (TB) 	to node [above right]	{$(\eta^{\mathbf{B}})^{\sharp}$} (FB);
\end{tikzpicture}
\end{center}
\caption{The free completion functor on a discrete $\Sigma$-algebra.}
\label{FFreeCompFunDisc}
\end{figure}


\begin{remark}
For every $S$-sorted set $X$, $\mathbf{T}_{\Sigma}(X)$ together with $\eta^{X}$ is the free completion of $\mathbf{D}_{\Sigma}(X)$.
\end{remark}

In the following proposition we state that the free completion of a partial $\Sigma$-algebra can also be characterized internally.

\begin{proposition}\label{PFreeCompInt}
Let $\mathbf{A}$ be a partial $\Sigma$-algebra and $\mathbf{B}$ a $\Sigma$-algebra such that $\mathbf{A}$ is a weak subalgebra of $\mathbf{B}$. Then $\mathbf{B}$ together with $\mathrm{in}^{\mathbf{A},\mathbf{B}}$, the canonical embedding of $\mathbf{A}$ into $\mathbf{B}$, is a free completion of $\mathbf{A}$ if, and only if, the following conditions hold:
\begin{enumerate}
\item[FC1.] For every $(\mathbf{s},s)\in S^{\star}\times S$, every $\sigma\in \Sigma_{\mathbf{s},s}$ and every $(b_{j})_{j\in \bb{\mathbf{s}}}\in B_{\mathbf{s}}$, if $\sigma^{\mathbf{B}}((b_{j})_{j\in \bb{\mathbf{s}}})\in A_{s}$, then $\sigma^{\mathbf{B}}((b_{j})_{j\in \bb{\mathbf{s}}}) = \sigma^{\mathbf{A}}((b_{j})_{j\in \bb{\mathbf{s}}})$ (thus, in particular, $(b_{j})_{j\in \bb{\mathbf{s}}}\in \mathrm{Dom}(\sigma^{\mathbf{A}})$).
\item[FC2.] For every $s\in S$, every $\sigma, \tau \in\Sigma_{\cdot,s}$, every $(b_{j})_{j\in \bb{\mathbf{s}}}\in B_{\mathbf{s}}$, where $\mathbf{s} = \mathrm{ar}(\sigma)$, and every $(c_{k})_{k\in \bb{\mathbf{t}}}\in B_{\mathbf{s}}$, where $\mathbf{t} = \mathrm{ar}(\tau)$, if 
$\sigma^{\mathbf{B}}((b_{j})_{j\in \bb{\mathbf{s}}}) = \tau^{\mathbf{B}}((c_{k})_{k\in \bb{\mathbf{t}}})\nin A_{s}$, then $\sigma = \tau$ and $(b_{j})_{j\in \bb{\mathbf{s}}} = (c_{k})_{k\in \bb{\mathbf{t}}}$ (i.e., ``outside of $\mathbf{A}$'' the second axiom of the notion of partial Dedekind-Peano is satisfied).
\item[FC3.] $\mathrm{Sg}_{\mathbf{B}}(A) = B$.
\end{enumerate}
\end{proposition}

\begin{remark}\label{PPAlgSubNorm}
The condition $\mathrm{FC1}$ in Proposition~\ref{PFreeCompInt} entails that the weak subalgebra $\mathbf{A}$ of $\mathbf{B}$ is even a relative subalgebra of $\mathbf{B}$. Relative subalgebras satisfying the condition $\mathrm{FC1}$ will be called \emph{normal} according to Schmidt~\cite{sch70}.
\end{remark}

\begin{remark}
There exists a natural transformation from the restriction of the functor $\mathbf{F}_{\Sigma}$ to $\mathsf{PAlg}(\Sigma)_{\mathrm{c}}$, the category of many-sorted partial $\Sigma$-algebras and closed homomorphisms, to the functor $(\cdot)^{\infty}$ from $\mathsf{PAlg}(\Sigma)_{\mathrm{c}}$ to $\mathsf{Alg}(\Sigma)$.
\end{remark}

\begin{figure}
$$\xymatrix@C=60pt@R=40pt{
    \mathbf{A}
    \ar@/_1pc/[ddr]_{f}
    \ar@/^1pc/[drr]^{\eta^{\mathbf{A}}}
    \ar[dr]|-{\mathrm{in}^{\mathbf{A},\mathbf{Sch}(f)}} \\
     & \mathbf{Sch}(f)
     \ar[d]_{f^{\mathrm{Sch}}}
     \ar[r]^-{\mathrm{in}^{\mathbf{Sch}(f)}}  &
     \mathbf{F}_{\Sigma}(\mathbf{A})
     \ar[d]^{f^{@}} \\
     &
     \mathbf{B}
     \ar[r]_{\eta^{\mathbf{B}}} &
     \mathbf{F}_{\Sigma}(\mathbf{B})
   }
$$
\caption{The Schmidt algebra relative to a partial algebra.}
\label{FSch}
\end{figure}

\begin{figure}
\begin{center}
\begin{tikzpicture}
[ACliment/.style={-{To [angle'=45, length=5.75pt, width=4pt, round]}
}]
\node[] (A) 	at 	(-3,1.5) 	[] 	{$\mathbf{D}_{\Sigma}(X)$};
\node[] (B) 	at 	(0,-2) 	[] 	{$\mathbf{B}$};
\node[]	(S)		at 	(0,0)	[]	{$\mathbf{Sch}(f)$};
\node[] (FA)	at	(4,0)	[]	{$\mathbf{T}_{\Sigma}(X)$};
\node[] (FB)	at	(4,-2)	[]	{$\mathbf{F}_{\Sigma}(\mathbf{B})$};
\draw[ACliment, bend right]  (A) 	to node [below left]	{$f$} (B);
\draw[ACliment, bend left]  (A) 	to node [above right]	{$\eta^{\mathbf{D}_{\Sigma}(X)}$} (FA);
\draw[ACliment]  (A) 	to node [above right]	{$\mathrm{in}^{\mathbf{D}_{\Sigma}(X),\mathbf{Sch}(f)}$} (S);
\draw[ACliment]  (B) 	to node [below]	{$\eta^{\mathbf{B}}$} (FB);
\draw[ACliment]  (S) 	to node [above]	{$\mathrm{in}^{\mathbf{Sch}(f)}$} (FA);
\draw[ACliment]  (S) 	to node [right]	{$f^{\mathrm{Sch}}$} (B);
\draw[ACliment]  (FA) 	to node [right]	{$f^{@}$} (FB);
\end{tikzpicture}
\end{center}
\caption{The Schmidt algebra for the discrete case.}
\label{FSchDisc}
\end{figure}

In connection with the functor $\mathbf{F}_{\Sigma}$ from $\mathsf{PAlg}(\Sigma)$ to $\mathsf{Alg}(\Sigma)$ we have the following generalized recursion theorem, which we will also call Schmidt's construction. The just mentioned theorem, as we will see later, is the basis of a model theoretic approach to the theory of many-sorted partial algebras. Moreover, relying on it, as we will see below, we associate to terms in a many-sorted algebra \emph{partial} term operations on the underlying many-sorted set of a many-sorted partial algebra, and we state the general homomorphism theorem of Schmidt. 


\begin{proposition}[The Schmidt construction]\label{PSch}
Let $f$ be a homomorphism from the partial $\Sigma$-algebra $\mathbf{A}$ to the partial $\Sigma$-algebra $\mathbf{B}$. Then there exists an $A$-generated relative subalgebra $\mathbf{Sch}(f)$ of $\mathbf{F}_{\Sigma}(\mathbf{A})$, the free completion of $\mathbf{A}$, and a closed homomorphism $f^{\mathrm{Sch}}\colon \mathbf{Sch}(f)\mor \mathbf{B}$ such that the  diagram in Figure~\ref{FSch} commutes, where $\mathrm{in}^{\mathbf{A},\mathbf{Sch}(f)}$ is the canonical inclusion of $\mathbf{A}$ into $\mathbf{Sch}(f)$. Moreover, $f^{\mathrm{Sch}}$ is the largest homomorphic extension of $f$ to an $A$-generated relative subalgebra of $\mathbf{F}_{\Sigma}(\mathbf{A})$, and it is the only closed one of this kind. In honour of J. Schmidt, who introduced these concepts in~\cite{sch70}, we will call $f^{\mathrm{Sch}}\colon \mathbf{Sch}(f)\mor \mathbf{B}$ the \emph{Schmidt closed} $\mathbf{A}$-\emph{initial extension} of $f$ and $\mathrm{Ker}(f^{\mathrm{Sch}})$, denoted by $\mathrm{SKer}(f)$, the \emph{Schmidt kernel} of $f$.
\end{proposition}

\begin{proof}
It suffices to take as $f^{\mathrm{Sch}}$ precisely $(\mathbf{Sch}(f),\Gamma\!_{f^{\mathrm{Sch}}},\mathbf{B})$, where $\Gamma\!_{f^{\mathrm{Sch}}}$, the underlying $S$-sorted function of $f^{\mathrm{Sch}}$, is $\mathrm{Sg}_{\mathbf{F}_{\Sigma}(\mathbf{A})\times \mathbf{B}}(\Gamma\!_{f})$ and $\mathbf{Sch}(f) = (f^{@})^{-1}[\mathbf{B}]$, where, with the customary abuse of notation, we have identified $B$ with $\eta^{\mathbf{B}}[B]$.
\end{proof}

%

\begin{remark}\label{RPAlgInit}
For a $\Sigma$-algebra $\mathbf{A}$, the $A$-generated relative subalgebras of 
$\mathbf{F}_{\Sigma}(\mathbf{A})$ are sometimes also called $\mathbf{A}$-\emph{initial segments} of $\mathbf{F}_{\Sigma}(\mathbf{A})$, and the closed congruences on $\mathbf{A}$-initial segments of $\mathbf{F}_{\Sigma}(\mathbf{A})$ are sometimes also called $\mathbf{A}$-\emph{initial congruences} of 
$\mathbf{F}_{\Sigma}(\mathbf{A})$.
\end{remark}
 
\begin{remark}
The canonical inclusion $\mathrm{in}^{\mathbf{A},\mathbf{Sch}(f)}$ of $\mathbf{A}$ into $\mathbf{Sch}(f)$ is an $\mathsf{Alg}(\Sigma)$-extendable epimorphism, where, we recall, a homomorphism $f\colon \mathbf{C}\mor \mathbf{D}$ in the category $\mathsf{PAlg}(\Sigma)$ is called $\mathsf{Alg}(\Sigma)$-\emph{extendable} if, for every $\Sigma$-algebra $\mathbf{E}$ in the category $\mathsf{Alg}(\Sigma)$ and every homomorphism $g\colon \mathbf{C}\mor \mathbf{E}$, there exists a homomorphism $h\colon \mathbf{D}\mor \mathbf{E}$ such that $h\circ f = g$.
%
\end{remark}

From Proposition~\ref{PSch} we obtain, for a homomorphism whose source is a discrete many-sorted partial algebra, the following corollary.

\begin{corollary}\label{CPAlgPartial}
Let $X$ be an $S$-sorted set, $\mathbf{B}$ a partial $\Sigma$-algebra, and $f$ an $S$-sorted mapping from $X$ to $B$ or, what is equivalent, a homomorphism from $\mathbf{D}_{\Sigma}(X)$ to $\mathbf{B}$. Then there exists an $X$-generated relative subalgebra $\mathbf{Sch}(f)$ of $\mathbf{T}_{\Sigma}(X)$, the free completion of $\mathbf{D}_{\Sigma}(X)$, and a closed homomorphism $$f^{\mathrm{Sch}}\colon \mathbf{Sch}(f)\mor \mathbf{B}$$ such that the diagram in Figure~\ref{FSchDisc} commutes, where $\mathrm{in}^{\mathbf{D}_{\Sigma}(X),\mathbf{Sch}(f)}$ is the canonical inclusion of $X$ into $\mathbf{Sch}(f)$. Moreover, $f^{\mathrm{Sch}}\colon \mathbf{Sch}(f)\mor \mathbf{B}$ is the largest homomorphic extension of $f$ to an $X$-generated relative subalgebra of $\mathbf{F}_{\Sigma}(X)$, and it is the only closed one of this kind.
\end{corollary}

Let us note that for $f\colon\mathbf{D}_{\Sigma}(X)\mor\mathbf{B}$ the diagram in Figure~\ref{FSch} becomes the diagram represented in Figure~\ref{FSchDisc}. Moreover, in this particular case, we have that $f^{@}=(\eta^{\mathbf{B}}\circ f)^{\sharp}$ and 
$\mathbf{Sch}(f)=((\eta^{\mathbf{B}}\circ f)^{\sharp})^{-1}[\eta^{\mathbf{B}}[\mathbf{B}]]$.

\begin{corollary}\label{CPAlgInitGen}
Let $\mathbf{B}$ be a partial $\Sigma$-algebra, $X$ a generating subset of $\mathbf{B}$, and 
$\mathrm{in}^{\mathbf{D}_{\Sigma}(X),\mathbf{B}}$ the canonical embedding of $\mathbf{D}_{\Sigma}(X)$ into $\mathbf{B}$. Then $\mathbf{B}$ is a closed homomorphic image of an $X$-generated relative subalgebra---namely of $\mathbf{Sch}(\mathrm{in}^{\mathbf{D}_{\Sigma}(X),\mathbf{B}})$---of $\mathbf{T}_{\Sigma}(X)$. (Therefore, for every $s\in S$, $\mathrm{card}(B_{s})\leq \mathrm{card}(\mathrm{T}_{\Sigma}(X)_{s})$.)
\end{corollary}

\begin{remark}\label{RSch}
We would like to highlight that in the proof of Proposition~\ref{PSch} as well as in Corollary~\ref{PPAlgPartial}, we have made use of the free completion both in the domain and in the codomain. One of the main reasons for the development of Climent and Cosme~\cite{CVCL23}---which aimed to show the functoriality of a generalization of the Schmidt construction---was to understand the fact that these free completions can be replaced by others with similar behavior. This led us, after defining a suitable category of completions and morphisms between them, to state and prove a generalization of the Schmidt construction. The reader can go to Proposition~\ref{PSch} to become aware of the generalisation of the Schmidt construction and to check the similarities between its proof and that of Proposition~\ref{PSch}.
\end{remark}

Also, on the basis of Proposition~\ref{PSch}, we associate to terms in a many-sorted algebra partial term operations, i.e., (generalized) partial operations, on the underlying many-sorted set of a many-sorted partial algebra. We point out that, as we shall see later on, partial term operations account for the semantics of the existence equations, the existentially conditioned existence equations, and the quasi-existence equations.

\begin{definition}
Let $X$ be an $S$-sorted set, $\mathbf{A}$ a partial $\Sigma$-algebra, $s\in S$ and $P\in \mathrm{T}_{\Sigma}(X)_{s}$. Then we denote by $P^{\mathbf{A}}$ the partial mapping from $A_{X} = \mathrm{Hom}(X,A)$ to $A_{s}$ defined as follows: 
$$
\mathrm{Dom}(P^{\mathbf{A}}) = \{a\in A_{X}\mid P\in \mathrm{Sch}(a)_{s}\}
$$ 
and, for every valuation $a\in \mathrm{Dom}(P^{\mathbf{A}})$ of $X$ in $A$, $P^{\mathbf{A}}(a) = a^{\mathrm{Sch}}_{s}(P)$. We will call $P^{\mathbf{A}}$ the \emph{partial term operation} on $\mathbf{A}$ induced by the term $P$. Moreover, we will say that $P$ is \emph{evaluable} in $\mathbf{A}$ with respect to a valuation $a\in A_{X}$ if and only if $P\in \mathrm{Sch}(a)_{s}$ (if and only if $a\in \mathrm{Dom}(P^{\mathbf{A}})$).
\end{definition}

\begin{remark}
Let $X$ be an $S$-sorted set of variables, $\mathbf{A}$ a partial $\Sigma$-algebra, $s\in S$ and $P\in \mathrm{T}_{\Sigma}(X)_{s}$.  Then we have that: If $P = x$, then $P^{\mathbf{A}}$ is the mapping $\mathrm{pr}^{A}_{X,s,x}$ from $A_{X}$ to $A_{s}$ that sends $a\in A_{X}$ to $\mathrm{pr}^{A}_{X,s,x}(a) = a_{s}(x)$; if $P = \sigma^{\mathbf{T}_{\Sigma}(X)}$ for some constant operation symbol $\sigma\in \Sigma_{\lambda,s}$ and $\mathrm{Dom}(\sigma^{\mathbf{A}})\neq\varnothing$, then $P^{\mathbf{A}}$ is the mapping from $A_{X}$ to $A_{s}$ that sends $a\in A_{X}$ to $\sigma^{\mathbf{A}}\in A_{s}$, while if $\mathrm{Dom}(\sigma^{\mathbf{A}}) = \varnothing$, then $\mathrm{Dom}(P^{\mathbf{A}}) = \varnothing$ and $P^{\mathbf{A}}$ is the unique mapping from $\varnothing$ to $A_{s}$; finally, if $P = \sigma^{\mathbf{T}_{\Sigma}(X)}((P_{i})_{i\in \bb{\mathbf{s}}})$, for $(\mathbf{s},s)\in S^{\star}\times S$ and $\sigma\in \Sigma_{\mathbf{s},s}$, then $P^{\mathbf{A}}$ is the partial mapping from $A_{X}$ to $A_{s}$, defined as follows:
$$
\textstyle
\mathrm{Dom}(P^{\mathbf{A}}) = \{a\in A_{X}\mid a\in \bigcap_{i\in \bb{\mathbf{s}}}\mathrm{Dom}(P_{i}^{\mathbf{A}}) \And (P^{\mathbf{A}}_{i}(a))_{i\in \bb{\mathbf{s}}}\in \mathrm{Dom}(\sigma^{\mathbf{A}})\},
$$ 
and, for every $a\in \mathrm{Dom}(P^{\mathbf{A}})$, $P^{\mathbf{A}}(a) = \sigma^{\mathbf{A}}((P^{\mathbf{A}}_{i}(a))_{i\in \bb{\mathbf{s}}})$. Thus $P^{\mathbf{A}}$ is the partial mapping
$\sigma^{\mathbf{A}}\circ \left<P^{\mathbf{A}}_{i}\right>_{i\in\bb{\mathbf{s}}}$ from $A_{X}$ to $A_{s}$, where the partial mapping $\left<P^{\mathbf{A}}_{i}\right>_{i\in\bb{\mathbf{s}}}$ from $A_{X}$ to $A_{\mathbf{s}}$ has the following universal property: It is the unique partial mapping from $A_{X}$ to $A_{\mathbf{s}}$ such that, for every $i\in \bb{\mathbf{s}}$, $\mathrm{pr}_{\mathbf{s},i}\circ  \left<P^{\mathbf{A}}_{i}\right>_{i\in\bb{\mathbf{s}}} \leq P^{\mathbf{A}}_{i}$, where $\mathrm{pr}_{\mathbf{s},i}$ is the $i$-th canonical projection from $A_{\mathbf{s}}$ to $A_{s_{i}}$, and, for every partial mapping $h$ from $A_{X}$ to $A_{s}$, if, for every $i\in \bb{\mathbf{s}}$, $\mathrm{pr}_{\mathbf{s},i}\circ h \leq P^{\mathbf{A}}_{i}$, then $h\leq \left<P^{\mathbf{A}}_{i}\right>_{i\in\bb{\mathbf{s}}}$ (we recall that for two partial mappings $u, v\colon X\dmor Y$, $u\leq v$ if, and only if, 
$\mathrm{Dom}(u)\subseteq \mathrm{Dom}(v)$ and, for every $x\in \mathrm{Dom}(u)$, $u(x) = v(x)$).

Let us also note that the second case is a particular instance of the last case. 
\end{remark}

Another very important tool provided by the generalized recursion theorem is the general homomorphism theorem of Schmidt. As a preparation we present the following characterization of epimorphisms.

\begin{lemma}
Let $f$ be a homomorphism from the partial $\Sigma$-algebra $\mathbf{A}$ to the partial $\Sigma$-algebra $\mathbf{B}$. Then $f$ is an epimorphism if and only if $f^{\mathrm{Sch}}\colon \mathbf{Sch}(f)\mor \mathbf{B}$ is surjective.
\end{lemma}

\begin{proof}
The reason why this is so is that 
$$
\mathrm{Sg}_{\mathbf{B}}(f[A]) = \mathrm{Sg}_{\mathbf{B}}(f^{\mathrm{Sch}}[A]) = f^{\mathrm{Sch}}[\mathrm{Sg}_{\mathbf{Sch}(f)}(A)] = f^{\mathrm{Sch}}[\mathrm{Sch}(f)].
$$
\end{proof}

\begin{proposition}
Let $f$ be a epimorphism from the partial $\Sigma$-algebra $\mathbf{A}$ to the partial $\Sigma$-algebra $\mathbf{B}$ and $g$ a homomorphism from $\mathbf{A}$ to the partial $\Sigma$-algebra $\mathbf{C}$. Then the following statements are equivalent:
\begin{enumerate}
\item[(i)] There exists a unique homomorphism $h$ from $\mathbf{B}$ to $\mathbf{C}$ such that $g = h\circ f$.
\item[(ii)] $\mathrm{SKer}(f)\subseteq \mathrm{SKer}(g)$ (as binary relations on the underlying $S$-sorted of
$\mathbf{F}_{\Sigma}(\mathbf{A})$).
\end{enumerate}
\end{proposition}

\begin{corollary}
Let $f$ be an epimorphism from the partial $\Sigma$-algebra $\mathbf{A}$ to the partial $\Sigma$-algebra $\mathbf{B}$. Then $\mathbf{B}\cong \mathbf{Sch}(f)/\mathrm{SKer}(f)$.
\end{corollary}

\section{Varieties of  many-sorted partial algebras}

The concept of existence equation is necessary for the syntactic description of certain classes of many-sorted  partial algebras, e.g., the categories, the $n$-categories, etc. 

\begin{definition}
Let $\Sigma$ be an $S$-sorted signature.
\begin{enumerate}
\item If $(X,s)\in\boldsymbol{\mathcal{U}}^{S}\times S$, then a
      $\Sigma$-\emph{term of type} $(X,s)$ is an
      element of $\mathrm{T}_{\Sigma}(X)_{s}$. If a $\Sigma$-term
      of type $(X,s)$ is such that $X$ is $S$-finite, i.e., for every $s\in S$, $X_{s}$ is finite, then we will
      call it $S$-\emph{finitary}, and if it is such that $X$ is
      finite, i.e., $\mathrm{supp}_{S}(X)$ is finite and, for every $s\in \mathrm{supp}_{S}(X)$, $X_{s}$ is 
      finite, then we will call it \emph{finitary}.

\item We denote by $\mathrm{Tm}(\Sigma)$ the $\boldsymbol{\mathcal{U}}^{S}\times S$-sorted
      set $(\mathrm{T}_{\Sigma}(X)_{s})_{(X,s)\in\boldsymbol{\mathcal{U}}^{S}\bprod S}$ 
      of all $\Sigma$-terms, by $\mathrm{Tm}_{S\text{-}\mathrm{f}}(\Sigma)$
      the $\boldsymbol{\mathcal{U}}^{S}_{S\text{-}\mathrm{f}}\times
      S$-sorted set 
      $(\mathrm{T}_{\Sigma}(X)_{s})_{(X,s)\in\boldsymbol{\mathcal{U}}^{S}_{S\text{-}\mathrm{f}}\bprod S}$ 
      of all $S$-finitary terms, where $\boldsymbol{\mathcal{U}}^{S}_{S\text{-}\mathrm{f}}$
      is the set of all $S$-sorted sets which are $S$-finite, and by
      $\mathrm{Tm}_{\mathrm{f}}(\Sigma)$ the $\boldsymbol{\mathcal{U}}^{S}_{\mathrm{f}}\times
      S$-sorted set 
      $(\mathrm{T}_{\Sigma}(X)_{s})_{(X,s)\in\boldsymbol{\mathcal{U}}^{S}_{\mathrm{f}}\bprod S}$ of all finitary terms, where $\boldsymbol{\mathcal{U}}^{S}_{\mathrm{f}}$ is
      the set of all $S$-sorted sets which are finite.
\end{enumerate}
\end{definition}

For $(X,s)\in\boldsymbol{\mathcal{U}}^{S}\times S$, we next define the notions of existence equation of type $(X,s)$, of $S$-finitary existence equation of type $(X,s)$ and of finitary existence equation of type $(X,s)$.

\begin{definition}
Let $\Sigma$ be an $S$-sorted signature.
\begin{enumerate}
\item If $(X,s)\in\boldsymbol{\mathcal{U}}^{S}\times S$, then an \emph{existence equation}
      (briefly $\mathrm{E}$-\emph{equation} or $\mathrm{Eeqt}$) of type $(X,s)$ is a pair $(P,Q)$, also denoted by
      $P\overset{\mathrm{e}}{=}Q$, where $P$ and $Q$ are $\Sigma$-terms of type $(X,s)$.
      If an $\mathrm{E}$-equation of type $(X,s)$ is such that $X$ is $S$-finite, then we will call it
      $S$\emph{-finitary}, and if it is such that $X$ is finite, then we will call it \emph{finitary}.

\item We denote by $\mathrm{Eeqt}(\Sigma)$ the $\boldsymbol{\mathcal{U}}^{S}\times S$-sorted set 
      $\mathrm{Tm}(\Sigma)^{2}$ of
      all $\mathrm{E}$-equations, by $\mathrm{Eeqt}_{S\text{-}\mathrm{f}}(\Sigma)$ the
      $\boldsymbol{\mathcal{U}}^{S}_{S\text{-}\mathrm{f}}\times S$-sorted set 
      $\mathrm{Tm}_{S\text{-}\mathrm{f}}(\Sigma)^{2}$ of all
      $S$-finitary $\mathrm{E}$-equations, by $\mathrm{Eeqt}_{\mathrm{f}}(\Sigma)$
      the $\boldsymbol{\mathcal{U}}^{S}_{\mathrm{f}}\times
      S$-sorted set $\mathrm{Tm}_{\mathrm{f}}(\Sigma)^{2}$ of all finitary $\mathrm{E}$-equations, and by $\mathrm{Eeqt}(\Sigma)_{X}$
      the $S$-sorted set $(\mathrm{T}_{\Sigma}(X)_{s}^{2})_{s\in S}$ of all $\mathrm{E}$-equations
      with variables in an $S$-sorted set $X$.
\item We will say that a family of $\mathrm{E}$-equations $\mathcal{E}\subseteq \mathrm{Eeqt}(\Sigma)$ is
      $S$-\emph{finitary} (resp. \emph{finitary}) if
      $\mathcal{E}\subseteq\mathrm{Eeqt}_{S\text{-} \mathrm{f}}(\Sigma)$
      (resp. $\mathcal{E}\subseteq\mathrm{Eeqt}_{\mathrm{f}}(\Sigma)$).
\end{enumerate}
\end{definition}

For $(X,s)\in\boldsymbol{\mathcal{U}}^{S}\times S$, we next define the ternary relation of satisfiability $\cdot \models^{\Sigma}_{X,s}\!\!\cdot\, [\cdot]$ between a partial $\Sigma$-algebra, an $\mathrm{E}$-equation of type $(X,s)$, and a valuation of $X$ in the underlying $S$-sorted set of the partial $\Sigma$-algebra, and the derived binary relation of validity $\cdot \models^{\Sigma}_{X,s}\!\!\cdot$ between a partial $\Sigma$-algebra and a $\mathrm{E}$-equation of type $(X,s)$. 

\begin{definition}
Let $\mathbf{A}$ be a partial $\Sigma$-algebra, $(X,s)\in\boldsymbol{\mathcal{U}}^{S}\times S$, $(P,Q)\in
\mathrm{Eeqt}(\Sigma)_{X,s}$, $\mathcal{E}\subseteq \mathrm{Eeqt}(\Sigma)$, and $a\in A_{X} = \mathrm{Hom}(X,A)$.
\begin{enumerate}
\item We will say that $\mathbf{A}$ \emph{satisfies the $\mathrm{E}$-equation} $(P,Q)$
      \emph{with respect to the valuation} $a$, denoted by 
      $$
      \mathbf{A}\models^{\Sigma}_{X,s}(P,Q) [a],
      $$
      if and only if $P\in \mathrm{Sch}(a)_{s}$, $Q\in \mathrm{Sch}(a)_{s}$,
      and $a^{\mathrm{Sch}}_{s}(P) = a^{\mathrm{Sch}}_{s}(Q)$.
\item We will say that $(P,Q)$ is \emph{valid} in $\mathbf{A}$, denoted by
      $$
      \mathbf{A}\models^{\Sigma}_{X,s}(P,Q),
      $$ 
      or by $\mathbf{A}\models(P,Q)$,
      when this is unlikely to cause confusion, if, for every $a\in A_{X}$,
      $\mathbf{A}\models^{\Sigma}_{X,s}(P,Q) [a]$.
\item We will say that $\mathcal{E}$ is \emph{valid} in $\mathbf{A}$, or that $\mathbf{A}$ \emph{satisfies} $\mathcal{E}$, denoted by
      $$
      \mathbf{A}\models^{\Sigma}\mathcal{E},
      $$ 
      or by $\mathbf{A}\models\mathcal{E}$,
      when this is unlikely to cause confusion, if, for every 
      $(X,s)\in\boldsymbol{\mathcal{U}}^{S}\times S$ and every $(P,Q)\in\mathcal{E}_{X,s}$, 
      $\mathbf{A}\models^{\Sigma}_{X,s}(P,Q)$.
\item We will denote by $\mathrm{Mod}_{\Sigma}(\mathcal{E})$, or by $\mathrm{PAlg}(\Sigma,\mathcal{E})$, the set $$
\{\mathbf{A}\in \mathrm{PAlg}(\Sigma)\mid \mathbf{A}\models^{\Sigma}\mathcal{E}\}
$$
of all partial $\Sigma$-algebras that satisfy $\mathcal{E}$.
\end{enumerate}
\end{definition}

\begin{remark}
Let $\mathbf{A}$ be a partial $\Sigma$-algebra, $(P,Q)\in\mathrm{Eeqt}(\Sigma)_{X,s}$, and $a\in A_{X}$.
If $P\in \mathrm{Sch}(a)_{s}$ and $Q\in \mathrm{Sch}(a)_{s}$, then, since $P^{\mathbf{A}}(a) = a^{\mathrm{Sch}}_{s}(P)$ and $Q^{\mathbf{A}}(a) = a^{\mathrm{Sch}}_{s}(P)$, we have that
$a^{\mathrm{Sch}}_{s}(P) = a^{\mathrm{Sch}}_{s}(Q)$ if and only if $P^{\mathbf{A}}(a) = Q^{\mathbf{A}}(a)$.
\end{remark}

\begin{definition}
Let $\mathbf{A}$ be a partial $\Sigma$-algebra and $\mathcal{E}\subseteq \mathrm{Eeqt}(\Sigma)$.
We will say that $\mathcal{E}$ is \emph{valid} in $\mathbf{A}$, denoted by
$\mathbf{A}\models^{\Sigma} \mathcal{E}$, if, for every $(X,s)\in
\boldsymbol{\mathcal{U}}^{S}\times S$ and every $(P,Q)\in \mathcal{E}_{X,s}$ we have that
$\mathbf{A}\models^{\Sigma}_{X,s} (P,Q)$.
\end{definition}

In some cases it is convenient to consider instead of the
$\boldsymbol{\mathcal{U}}^{S}\times S$-sorted set $\mathrm{Eeqt}(\Sigma)$ the set
$\coprod\mathrm{Eeqt}(\Sigma) = \bigcup_{(X,s)\in\boldsymbol{\mathcal{U}}^{S}\times S}(\mathrm{T}_{\Sigma}(X)_{s}^{2}\times\{(X,s)\})$. Since
$$\mathbf{\Sub}(\mathrm{Eeqt}(\Sigma))\iso\mathbf{\Sub}(\coprod\mathrm{Eeqt}(\Sigma)),$$ it follows that
every family $\mathcal{E}\subseteq\mathrm{Eeqt}(\Sigma)$ has univocally associated a subset of $\coprod\mathrm{Eeqt}(\Sigma)$ and reciprocally. In the same way, sometimes, it is convenient to consider instead of the $S$-sorted set $\mathrm{Eeqt}(\Sigma)_{X}$ the set $\coprod\mathrm{Eeqt}(\Sigma)_{X}$.

\begin{remark}
$\mathrm{Sub}(\mathrm{Eeqt}(\Sigma))$ and $\prod_{X\in
\boldsymbol{\mathcal{U}}^{S}}\mathrm{Sub}(\coprod\mathrm{Eeqt}(\Sigma)_{X})$ are naturally isomorphic. In fact, it suffices to consider the mapping that assigns to $\mathcal{E}\subseteq \mathrm{Eeqt}(\Sigma)$ precisely 
$(\coprod_{s\in S}\mathcal{E}_{X,s})_{X\in \boldsymbol{\mathcal{U}}^{S}}$. 
\end{remark}

For $(X,(\mathbf{s},s))$ in $\boldsymbol{\mathcal{U}}^{S}\times (S^{\star}\times S)$, we next define the notions of existentially conditioned existence equation of type $(X,(\mathbf{s},s))$, of $S$-finitary existentially conditioned existence equation of type $(X,(\mathbf{s},s))$ and of finitary existentially conditioned existence equation of type $(X,(\mathbf{s},s))$.

\begin{definition}
Let $\Sigma$ be an $S$-sorted signature.
\begin{enumerate}
\item If $(X,(\mathbf{s},s))\in\boldsymbol{\mathcal{U}}^{S}\times (S^{\star}\times S)$, then an \emph{existentially conditioned existence equation} (briefly $\mathrm{ECE}$-\emph{equation} or $\mathrm{ECEeqt}$) of type $(X,(\mathbf{s},s))$ is a pair
$$
\left(\left(
P_{j},P_{j}
\right)_{j\in \bb{\mathbf{s}}},(P,Q)\right),
$$
also denoted by
$$
\textstyle
\bigwedge_{j\in\bb{\mathbf{s}}}(P_{j},P_{j})\to (P,Q)\, \text{ or by }\,
\bigwedge_{j\in\bb{\mathbf{s}}}(P_{j}\overset{\mathrm{e}}{=}P_{j})\to (P\overset{\mathrm{e}}{=}Q),
$$
where, for every $j\in \bb{\mathbf{s}}$, $P_{j}\in\mathrm{T}_{\Sigma}(X)_{s_{j}}$ and $(P,Q)\in \mathrm{T}_{\Sigma}(X)^{2}_{s}$.
The pairs $(P_{j},P_{j})$, for $j\in \bb{\mathbf{s}}$, are called the \emph{premisses} and $(P,Q)$ the \emph{conclusion} of the $\mathrm{ECE}$-equation. If an $\mathrm{ECE}$-equation of type $(X,(\mathbf{s}, s))$ is such that $X$ is $S$-finite, then we will call it $S$-\emph{finitary}, and if it is such that $X$ is finite, then we will call it \emph{finitary}.
\item We denote by $\mathrm{ECEeqt}(\Sigma)_{X}$ the $S^{\star}\times S$-sorted set of all 
      $\mathrm{ECE}$-equations with
      variables in $X$, by $\mathrm{ECEeqt}(\Sigma)$
      the $\boldsymbol{\mathcal{U}}^{S}\times (S^{\star}\times S)$-sorted set of
      all $\mathrm{ECE}$-equations, by $\mathrm{ECEeqt}_{S\text{-}\mathrm{f}}(\Sigma)$ the
      $\boldsymbol{\mathcal{U}}^{S}_{S\text{-}\mathrm{f}}\times (S^{\star}\times S)$-sorted set of all
      $S$-finitary $\mathrm{ECE}$-equations, and by $\mathrm{ECEeqt}_{\mathrm{f}}(\Sigma)$
      the $\boldsymbol{\mathcal{U}}^{S}_{\mathrm{f}}\times (S^{\star}\times S)$-sorted set of all finitary
      $\mathrm{ECE}$-equations.
\item We will say that a family of $\mathrm{ECE}$-equations $\mathcal{E}\subseteq \mathrm{ECEeqt}(\Sigma)$ is
      $S$-\emph{finitary} (resp. \emph{finitary}) if
      $\mathcal{E}\subseteq\mathrm{ECEeqt}_{S\text{-} \mathrm{f}}(\Sigma)$
      (resp. $\mathcal{E}\subseteq\mathrm{ECEeqt}_{\mathrm{f}}(\Sigma)$).
\end{enumerate}
\end{definition}

In some cases it is convenient to consider instead of the
$\boldsymbol{\mathcal{U}}^{S}\times (S^{\star}\times S)$-sorted set $\mathrm{ECEeqt}(\Sigma)$ the set $\coprod\mathrm{ECEeqt}(\Sigma)$, which is
$$
\textstyle
\bigcup_{(X,(\mathbf{s},s))\in\boldsymbol{\mathcal{U}}^{S}\times (S^{\star}\times S)}
\left(\left(\left(\prod_{j\in\bb{\mathbf{s}}}\Delta_{\mathrm{T}_{\Sigma}(X)_{s_{j}}}\right)\times\mathrm{T}_{\Sigma}(X)_{s}^{2}\right)\times\{(X,(\bb{\mathbf{s}},s))\}\right).
$$
Since
$\mathbf{\Sub}(\mathrm{ECEeqt}(\Sigma))\iso\mathbf{\Sub}(\coprod\mathrm{ECEeqt}(\Sigma))$, it follows that
every family $\mathcal{E}\subseteq\mathrm{ECEeqt}(\Sigma)$ has univocally associated a subset of $\coprod\mathrm{ECEeqt}(\Sigma)$ and reciprocally. In the same way, sometimes, it is convenient to consider instead of the $(S^{\star}\times S)$-sorted set $\mathrm{ECEeqt}(\Sigma)_{X}$ the set $\coprod\mathrm{ECEeqt}(\Sigma)_{X}$ which is
$$
\textstyle
\bigcup_{(\mathbf{s},s)\in S^{\star}\times S}
\left(\left(\left(
\prod_{j\in\bb{\mathbf{s}}}\Delta_{\mathrm{T}_{\Sigma}(X)_{s_{j}}}
\right)\times \mathrm{T}_{\Sigma}(X)^{2}_{s}
\right)\times \{(\mathbf{s},s)\}\right),
$$
which, in its turn, is isomorphic to
$$
\textstyle
\left(
\coprod_{s\in S}\Delta_{\mathrm{T}_{\Sigma}(X)_{s}}
\right)^{\star}\times \coprod_{s\in S}\mathrm{T}_{\Sigma}(X)_{s}^{2}.
$$

For $(X,(\mathbf{s},s))\in\boldsymbol{\mathcal{U}}^{S}\times (S^{\star}\times S)$, we next define the ternary relation of satisfiability $\cdot \models^{\Sigma}_{X,\mathbf{s},s}\!\!\cdot\, [\cdot]$ between a partial $\Sigma$-algebra, an $\mathrm{ECE}$-\emph{equation} of type $(X,(\mathbf{s},s))$, and a valuation of $X$ in the underlying $S$-sorted set of the partial $\Sigma$-algebra, and the derived binary relation of validity $\cdot \models^{\Sigma}_{X,\mathbf{s},s}\!\!\cdot$ between a partial $\Sigma$-algebra anf a $\mathrm{ECE}$-equation of type $(X,(\mathbf{s},s))$. 

\begin{definition}
Let $\mathbf{A}$ be a partial $\Sigma$-algebra, $(X,(\mathbf{s},s))\in\boldsymbol{\mathcal{U}}^{S}\times (S^{\star}\times S)$, $\bigwedge_{j\in\bb{\mathbf{s}}}(P_{j},P_{j})\to (P,Q)\in
\mathrm{ECEeqt}(\Sigma)_{(X,(\mathbf{s},s))}$, $\mathcal{E}\subseteq \mathrm{ECEeqt}(\Sigma)$, and $a\in A_{X} = \mathrm{Hom}(X,A)$.
\begin{enumerate}
\item We will say that $\mathbf{A}$ \emph{satisfies the $\mathrm{ECE}$-equation} 
      $\bigwedge_{j\in\bb{\mathbf{s}}}(P_{j},P_{j})\to (P,Q)$
      \emph{with respect to the valuation} $a$ , denoted by 
      $$
      \mathbf{A}\models^{\Sigma}_{X,\mathbf{s},s}\bigwedge_{j\in\bb{\mathbf{s}}}(P_{j},P_{j})\to (P,Q) [a],
      $$
      if and only if, if, for every $j\in \bb{\mathbf{s}}$, 
      $P_{j}\in \mathrm{Sch}(a)_{s_{j}}$ (i.e., $P^{\mathbf{A}}_{j}$ is defined at $a$), then 
      $P\in \mathrm{Sch}(a)_{s}$, $Q\in \mathrm{Sch}(a)_{s}$,
      and $a^{\mathrm{Sch}}_{s}(P) = a^{\mathrm{Sch}}_{s}(Q)$. In other words, if, for every 
      $j\in \bb{\mathbf{s}}$, $\mathbf{A}\models^{\Sigma}_{X,s_{j}} (P_{j},P_{j})[a]$, then  
      $\mathbf{A}\models^{\Sigma}_{X,s} (P,Q)[a]$.
\item We will say that $\bigwedge_{j\in\bb{\mathbf{s}}}(P_{j},P_{j})\to (P,Q)$ is \emph{valid} in 
      $\mathbf{A}$, denoted by
      $$
      \mathbf{A}\models^{\Sigma}_{X,\mathbf{s},s}\bigwedge_{j\in\bb{\mathbf{s}}}(P_{j},P_{j})\to (P,Q),
      $$ 
      or by $\mathbf{A}\models\bigwedge_{j\in\bb{\mathbf{s}}}(P_{j},P_{j})\to (P,Q)$,
      when this is unlikely to cause confusion, if, for every $a\in A_{X}$,
      $\mathbf{A}\models^{\Sigma}_{X,\mathbf{s},s}\bigwedge_{j\in\bb{\mathbf{s}}}(P_{j},P_{j})\to (P,Q) [a]$.
\item We will say that $\mathcal{E}$ is \emph{valid} in $\mathbf{A}$, or that $\mathbf{A}$ \emph{satisfies} $\mathcal{E}$, denoted by
      $$
      \mathbf{A}\models^{\Sigma}\mathcal{E},
      $$
      or by $\mathbf{A}\models\mathcal{E}$,
      when this is unlikely to cause confusion, if, for every 
      $(X,(\mathbf{s},s))\in\boldsymbol{\mathcal{U}}^{S}\times (S^{\star}\times S)$ and every 
      $\bigwedge_{j\in\bb{\mathbf{s}}}(P_{j},P_{j})\to (P,Q)\in
      \mathcal{E}_{(X,(\mathbf{s},s))}$, 
      $$
      \mathbf{A}\models^{\Sigma}_{X,\mathbf{s},s}\bigwedge_{j\in\bb{\mathbf{s}}}(P_{j},P_{j})\to (P,Q).
      $$
\item We will denote by $\mathrm{Mod}_{\Sigma}(\mathcal{E})$, or by $\mathrm{PAlg}(\Sigma,\mathcal{E})$, the set 
$$
\{\mathbf{A}\in \mathrm{PAlg}(\Sigma)\mid \mathbf{A}\models^{\Sigma}\mathcal{E}\}
$$
of all partial $\Sigma$-algebras that satisfy $\mathcal{E}$.
\end{enumerate}
\end{definition}

\begin{definition}
Let $\mathbf{A}$ be a partial $\Sigma$-algebra and $\mathcal{E}\subseteq \mathrm{ECEeqt}(\Sigma)$.
We will say that $\mathcal{E}$ is \emph{valid} in $\mathbf{A}$, denoted by
$\mathbf{A}\models^{\Sigma} \mathcal{E}$, if, for every $(X,(\mathbf{s},s))\in
\boldsymbol{\mathcal{U}}^{S}\times (S^{\star}\times S)$ and every $\bigwedge_{j\in\bb{\mathbf{s}}}(P_{j},P_{j})\to (P,Q)\in \mathcal{E}_{(X,(\mathbf{s},s))}$ we have that
$\mathbf{A}\models^{\Sigma}_{X,\mathbf{s},s} \bigwedge_{j\in\bb{\mathbf{s}}}(P_{j},P_{j})\to (P,Q)$.
\end{definition}

For $(X,(\mathbf{s},s))$ in $\boldsymbol{\mathcal{U}}^{S}\times (S^{\star}\times S)$, we next define the notions of  quasi-existence equation of type $(X,(\mathbf{s},s))$, of $S$-finitary  quasi-existence equation of type $(X,(\mathbf{s},s))$ and of finitary  quasi-existence equation of type $(X,(\mathbf{s},s))$.

\begin{definition}\label{DQEEq}
Let $\Sigma$ be an $S$-sorted signature.
\begin{enumerate}
\item If $(X,(\mathbf{s},s))\in\boldsymbol{\mathcal{U}}^{S}\times (S^{\star}\times S)$, then a 
\emph{quasi-existence equation} (briefly $\mathrm{QE}$-\emph{equation} or $\mathrm{QEeqt}$) of type $(X,(\mathbf{s},s))$ is a pair
$$
\left(\left(
P_{j},Q_{j}
\right)_{j\in \bb{\mathbf{s}}},(P,Q)\right),
$$
also denoted by
$$
\textstyle
\bigwedge_{j\in\bb{\mathbf{s}}}(P_{j},Q_{j})\to (P,Q)\, \text{ or by }\,
\bigwedge_{j\in\bb{\mathbf{s}}}(P_{j}\overset{\mathrm{e}}{=}Q_{j})\to (P\overset{\mathrm{e}}{=}Q),
$$
where, for every $j\in \bb{\mathbf{s}}$, $(P_{j},Q_{j})\in\mathrm{T}_{\Sigma}(X)^{2}_{s_{j}}$ and $(P,Q)\in \mathrm{T}_{\Sigma}(X)^{2}_{s}$.
The pairs $(P_{j},Q_{j})$, for $j\in \bb{\mathbf{s}}$, are called the \emph{premisses} and $(P,Q)$ the \emph{conclusion} of the $\mathrm{QE}$-equation. If a $\mathrm{QE}$-equation of type $(X,(\mathbf{s}, s))$ is such that $X$ is $S$-finite, then we will call it $S$-\emph{finitary}, and if it is such that $X$ is finite, then we will call it \emph{finitary}.
\item We denote by $\mathrm{QEeqt}(\Sigma)_{X}$ the $S^{\star}\times S$-sorted set of all 
      $\mathrm{QE}$-equations with
      variables in $X$, by $\mathrm{QEeqt}(\Sigma)$
      the $\boldsymbol{\mathcal{U}}^{S}\times (S^{\star}\times S)$-sorted set of
      all $\mathrm{QE}$-equations, by $\mathrm{QEeqt}_{S\text{-}\mathrm{f}}(\Sigma)$ the
      $\boldsymbol{\mathcal{U}}^{S}_{S\text{-}\mathrm{f}}\times (S^{\star}\times S)$-sorted set of all
      $S$-finitary $\mathrm{QE}$-equations, and by $\mathrm{QEeqt}_{\mathrm{f}}(\Sigma)$
      the $\boldsymbol{\mathcal{U}}^{S}_{\mathrm{f}}\times (S^{\star}\times S)$-sorted set of all finitary
      $\mathrm{QE}$-equations.
\item We will say that a family of $\mathrm{QE}$-equations $\mathcal{E}\subseteq \mathrm{QEeqt}(\Sigma)$ is
      $S$-\emph{finitary} (resp. \emph{finitary}) if
      $\mathcal{E}\subseteq\mathrm{QEeqt}_{S\text{-} \mathrm{f}}(\Sigma)$
      (resp. $\mathcal{E}\subseteq\mathrm{QEeqt}_{\mathrm{f}}(\Sigma)$).
\end{enumerate}
\end{definition}

In some cases it is convenient to consider instead of the
$\boldsymbol{\mathcal{U}}^{S}\times (S^{\star}\times S)$-sorted set $\mathrm{QEeqt}(\Sigma)$ the set $\coprod\mathrm{QEeqt}(\Sigma)$, which is
$$
\textstyle
\bigcup_{(X,(\mathbf{s},s))\in\boldsymbol{\mathcal{U}}^{S}\times (S^{\star}\times S)}
\left(\left(\left(\prod_{j\in\bb{\mathbf{s}}}\mathrm{T}_{\Sigma}(X)_{s_{j}}\right)\times\mathrm{T}_{\Sigma}(X)_{s}^{2}\right)\times\{(X,(\bb{\mathbf{s}},s))\}\right).
$$
Since
$\mathbf{\Sub}(\mathrm{QEeqt}(\Sigma))\iso\mathbf{\Sub}(\coprod\mathrm{QEeqt}(\Sigma))$, it follows that
every family $\mathcal{E}\subseteq\mathrm{QEeqt}(\Sigma)$ has univocally associated a subset of $\coprod\mathrm{QEeqt}(\Sigma)$ and reciprocally. In the same way, sometimes, it is convenient to consider instead of the $(S^{\star}\times S)$-sorted set $\mathrm{QEeqt}(\Sigma)_{X}$ the set $\coprod\mathrm{QEeqt}(\Sigma)_{X}$ which is
$$
\textstyle
\bigcup_{(\mathbf{s},s)\in S^{\star}\times S}
\left(\left(\left(
\prod_{j\in\bb{\mathbf{s}}}\mathrm{T}_{\Sigma}(X)_{s_{j}}
\right)\times \mathrm{T}_{\Sigma}(X)^{2}_{s}
\right)\times \{(\mathbf{s},s)\}\right),
$$
which, in its turn, is isomorphic to
$$
\textstyle
\left(
\coprod_{s\in S}\mathrm{T}_{\Sigma}(X)^{2}_{s}
\right)^{\star}\times \coprod_{s\in S}\mathrm{T}_{\Sigma}(X)_{s}^{2}.
$$

Since the semantics of the $\mathrm{QE}$-equations follows the same pattern that the semantics of the $\mathrm{ECE}$-equations, we leave the details to the reader. Let us just say that a partial $\Sigma$-algebra $\mathbf{A}$ satisfies a $\mathrm{QE}$-equation 
      $\bigwedge_{j\in\bb{\mathbf{s}}}(P_{j},Q_{j})\to (P,Q)$
      with respect to a valuation $a$, denoted by 
      $$
      \mathbf{A}\models^{\Sigma}_{X,\mathbf{s},s}\bigwedge_{j\in\bb{\mathbf{s}}}(P_{j},Q_{j})\to (P,Q) [a],
      $$
      if and only if, if, for every $j\in \bb{\mathbf{s}}$, 
      $P_{j}\in \mathrm{Sch}(a)_{s_{j}}$, $Q_{j}\in \mathrm{Sch}(a)_{s_{j}}$,
      and $a^{\mathrm{Sch}}_{s_{j}}(P_{j}) = a^{\mathrm{Sch}}_{s_{j}}(Q_{j})$, then 
      $P\in \mathrm{Sch}(a)_{s}$, $Q\in \mathrm{Sch}(a)_{s}$,
      and $a^{\mathrm{Sch}}_{s}(P) = a^{\mathrm{Sch}}_{s}(Q)$. In other words, if, for every 
      $j\in \bb{\mathbf{s}}$, $\mathbf{A}\models^{\Sigma}_{X,s_{j}} (P_{j},Q_{j})[a]$, then  
      $\mathbf{A}\models^{\Sigma}_{X,s} (P,Q)[a]$.


We next define, respecting the terminology used by Burmeister in~\cite{pb86}, the notions of finitary $\mathrm{E}$-, finitary $\mathrm{ECE}$-, and finitary $\mathrm{QE}$-varieties---which will be the type of varieties to which we will restrict ourselves in this work---and, for this purpose, we, on the one hand, chose a standard $S$-countably infinite $S$-sorted set of variables $V^{S} = (\{v^{s}_{n}\mid n\in\mathbb{N}\})_{s\in S}$ (one could take as $V^{S}$ the $S$-sorted set $(\{s\}\times\mathbb{N})_{s\in S}$) and, on the other hand, assume that all involved $S$-sorted sets of variables $X$ are finite subsets of $V^{S}$, i.e., $\mathrm{supp}_{S}(X)$ is finite and, for every $s\in \mathrm{supp}_{S}(X)$, $X_{s}$ is finite (therefore, $X$ will be isomorphic to $\downarrow\!\! w$, for some $w\in S^{\star}$).

\begin{definition}
We will say that a set $\boldsymbol{\mathcal{K}}$ of partial $\Sigma$-algebras is a \emph{finitary} $\mathrm{E}$-, \emph{finitary} $\mathrm{ECE}$-, or \emph{finitary} $\mathrm{QE}$-\emph{variety}, respectively, if there exists a set $\mathcal{E}$ of finitary $\mathrm{E}$-, $\mathrm{ECE}$-, or $\mathrm{QE}$-equations, respectively, such that $\boldsymbol{\mathcal{K}} = \mathrm{Mod}_{\Sigma}(\mathcal{E})$.
\end{definition}

\begin{remark}
The above finitary varieties are also referred to as \emph{finitary} $\mathrm{E}$-, \emph{finitary} $\mathrm{ECE}$-, or \emph{finitary} $\mathrm{QE}$-\emph{equational class}, respectively. Moreover, for the characterization of the aforementioned varieties by closure properties with respect to algebraic constructions we refer the reader to~\cite{pb86}. 
\end{remark}

\begin{remark}
Categories (in a fixed Grothendieck universe) constitute a finitary $\mathrm{ECE}$-variety and not only a finitary $\mathrm{QE}$-variety (see~\cite{pb02}).
\end{remark}

We next state the fundamental theorem on the existence of universal solutions for sets of partial algebras closed under subalgebras and products. But before doing that, for partial $\Sigma$-algebras $\mathbf{A}$ and $\mathbf{B}$, we make the convention that $f\colon \mathbf{B}\epi\mathbf{A}$ means that $f$ is a surjective homomorphism from $\mathbf{B}$ to $\mathbf{A}$ and that $f\colon \mathbf{A}\mon\mathbf{B}$ means that $f$ is a a closed and injective homomorphism from $\mathbf{A}$ to $\mathbf{B}$.

\begin{definition}
Let $\mathcal{K}$ be a set of partial $\Sigma$-algebras. Then
\begin{enumerate}
\item $\mathrm{I}(\mathcal{K}) = \{\mathbf{A}\in \mathrm{Alg}(\Sigma)\mid \exists\,\mathbf{B}\in \mathcal{K}\,(\mathbf{A}\cong\mathbf{B})\}$, i.e., $\mathrm{I}(\mathcal{K})$ is the set of all isomorphic copies of partial $\Sigma$-algebras in $\mathcal{K}$.

\item $\mathrm{H}(\mathcal{K}) = \{\mathbf{A}\in \mathrm{Alg}(\Sigma)\mid \exists\,\mathbf{B}\in \mathcal{K}\,\exists\,f\colon \mathbf{B}\epi\mathbf{A}\}$, i.e., $\mathrm{H}(\mathcal{K})$ is the set of all homomorphic images of partial $\Sigma$-algebras in $\mathcal{K}$.

\item $\mathrm{S}(\mathcal{K}) = \{\mathbf{A}\in \mathrm{Alg}(\Sigma)\mid \exists\,\mathbf{B}\in \mathcal{K}\,\exists\,f\colon \mathbf{A}\mon\mathbf{B}\}$, i.e., $\mathrm{S}(\mathcal{K})$ is the set of all isomorphic copies of subalgebras of partial $\Sigma$-algebras in $\mathcal{K}$.

\item $\mathrm{P}(\mathcal{K}) = \{\mathbf{A}\in \mathrm{Alg}(\Sigma)\mid \exists\, I\in \boldsymbol{\mathcal{U}}\, \exists\, (\mathbf{A}^{i})_{i\in I}\in \mathcal{K}^{I}\,(\mathbf{A}\cong \prod_{i\in I}\mathbf{A}^{i})\}$, i.e., $\mathrm{P}(\mathcal{K})$ is the set of all partial $\Sigma$-algebras which are  isomorphic to products of families, indexed by some $\boldsymbol{\mathcal{U}}$-small set, of partial $\Sigma$-algebras in $\mathcal{K}$.
\item We will say that $\mathcal{K}$ is (1) \emph{abstract} if $\mathrm{I}(\mathcal{K})\subseteq \mathcal{K}$, (2) \emph{closed under homomorphic images} if $\mathrm{H}(\mathcal{K})\subseteq \mathcal{K}$, (3) \emph{closed under subalgebras} if $\mathrm{S}(\mathcal{K})\subseteq \mathcal{K}$, and (4) \emph{closed under products} if $\mathrm{P}(\mathcal{K})\subseteq \mathcal{K}$. 
\end{enumerate}
\end{definition} 

\begin{lemma}
Let $\mathbf{A}$ be a partial $\Sigma$-algebra. Then:
\begin{enumerate}
\item The set of all $A$-generated relative subalgebras of $\mathbf{F}_{\Sigma}(\mathbf{A})$ is a closure system on $\mathrm{F}_{\Sigma}(\mathbf{A})$ with smallest element $\mathbf{A}$ and greatest element $\mathbf{F}_{\Sigma}(\mathbf{A})$.
\item The set of all closed congruences on $A$-generated relative subalgebras of $\mathbf{F}_{\Sigma}(\mathbf{A})$ is a closure system on $\mathrm{F}_{\Sigma}(\mathbf{A})^{2}$ with smallest element $\Delta_{\mathbf{A}}$ and greatest element $\nabla_{{\mathbf{F}_{\Sigma}(\mathbf{A})}}$.
\item For every nonempty family $(\Phi^{i})_{i\in I}$ of closed congruences on $A$-generated relative subalgebras of $\mathbf{F}_{\Sigma}(\mathbf{A})$, we have that
$$
\textstyle
\mathrm{Fld}(\bigcap_{i\in I}\Phi^{i}) = \bigcap_{i\in I}\mathrm{Fld}(\Phi^{i}).
$$
\item The set of all closed congruences on $A$-generated relative subalgebras of $\mathbf{F}_{\Sigma}(\mathbf{A})$, ordered by inclusion, is upwards directed.
\end{enumerate}
\end{lemma}

\begin{corollary}
Let $\mathbf{A}$ be a partial $\Sigma$-algebra, $(\mathbf{B}^{i})_{i\in I}$ an $I$-indexed family of partial $\Sigma$-algebras, $(f^{i})_{i\in I}\in \prod_{i\in I}\mathrm{Hom}(\mathbf{A},\mathbf{B}^{i})$ and 
$\left< f^{i}\right>_{i\in I}$ the unique homomorphism from $\mathbf{A}$ to 
$\prod_{i\in I}\mathbf{B}^{i}$ such that, for every $i\in I$, 
$\mathrm{pr}^{i}\circ \left< f^{i}\right>_{i\in I} = f^{i}$. Then
\begin{enumerate}
\item $\mathbf{Sch}(\left< f^{i}\right>_{i\in I}) = \bigcap_{i\in I}\mathbf{Sch}(f^{i})$.
\item $\mathrm{SKer}(\left< f^{i}\right>_{i\in I}) = \bigcap_{i\in I}\mathrm{SKer}(f^{i})$.
\end{enumerate}
\end{corollary}

\begin{lemma}\label{LCardinal}
For every $S$-cardinal $(\mathfrak{n}_{s})_{s\in S}$, there exists a cardinal $\mathfrak{m}$ such that, for every partial $\Sigma$-algebra $\mathbf{B}$ and every $S$-sorted set $Y$, if $Y\subseteq B$, $\mathrm{Sg}_{\mathbf{B}}(Y) = B$ and $(\mathrm{card}(Y_{s}))_{s\in S}\leq (\mathfrak{n}_{s})_{s\in S}$, then $(\mathrm{card}(B_{s}))_{s\in S}\leq (\mathfrak{m})_{s\in S}$, i.e., for every $s\in S$, $\mathrm{card}(B_{s})\leq \mathfrak{m}$.
\end{lemma}

\begin{proof}
It sufices to take $\mathfrak{m} = \mathrm{max}(\{\mathrm{card}(\Sigma),\aleph_{0}\}\cup\{\mathfrak{n}_{s}\mid s\in S\})$ (by virtue of Corollary~\ref{CPAlgInitGen} and the recursive definition of the free $\Sigma$-algebra on $Y$).
\end{proof}


\begin{theorem}\label{TFreeAdj}
Let $\mathcal{K}$ be an abstract set of partial $\Sigma$-algebras closed under subalgebras and products. Then the inclusion functor $\mathrm{In}_{\boldsymbol{\mathcal{K}}}$ from $\boldsymbol{\mathcal{K}}$, the full subcategory of $\mathsf{PAlg}(\Sigma)$ determined by $\mathcal{K}$, to $\mathsf{PAlg}(\Sigma)$ has a left adjoint $\mathbf{T}_{\boldsymbol{\mathcal{K}}}$.
\end{theorem}

\begin{proof}

The proof will consists in showing that for every partial $\Sigma$-algebra $\mathbf{A}$ there exists a universal morphism from $\mathbf{A}$ to $\mathrm{In}_{\boldsymbol{\mathcal{K}}}$.

\textsf{First proof (via Schmidt)}.\hfill

Let $\mathbf{A}$ be a partial $\Sigma$-algebra and $\equiv^{\mathcal{K}}_{\mathbf{A}}$ be
$\bigcap_{f\in\bigcup_{\mathbf{K}\in\mathcal{K}}\mathrm{Hom}(\mathbf{A},\mathbf{K})}\mathrm{SKer}(f)$ (recall that $\mathrm{SKer}(f) = \mathrm{Ker}(f^{\mathrm{Sch}})$ and that it is a closed congruence on the $A$-generated relative subalgebras $\mathbf{Sch}(f)$ of $\mathbf{F}_{\Sigma}(\mathbf{A})$). For every
$f\in\bigcup_{\mathbf{K}\in\mathcal{K}}\mathrm{Hom}(\mathbf{A},\mathbf{K})$, $\mathrm{Fld}(\mathrm{SKer}(f))$, the field of $\mathrm{SKer}(f)$, is $\mathrm{Sch}(f)$ (recall that, for an $S$-sorted equivalence relation $\Phi$ on an $S$-sorted set $Z$, since $\Phi$ is symmetric, $\mathrm{Fld}(\Phi) = \mathrm{Dom}(\Phi) = \mathrm{Im}(\Phi)$, and, since $\Phi$ is reflexive, $\mathrm{Fld}(\Phi) = Z$). Moreover
$$
\textstyle
\mathrm{Fld}
\left(
\bigcap_{f\in\bigcup_{\mathbf{K}\in\mathcal{K}}\mathrm{Hom}(\mathbf{A},\mathbf{K})}\mathrm{SKer}(f)
\right) =
\bigcap_{f\in\bigcup_{\mathbf{K}\in\mathcal{K}}\mathrm{Hom}(\mathbf{A},\mathbf{K})}\mathrm{Fld}(\mathrm{SKer}(f)).
$$
Therefore
$$
\textstyle
\mathrm{Fld}\left(
\bigcap_{f\in\bigcup_{\mathbf{K}\in\mathcal{K}}\mathrm{Hom}(\mathbf{A},\mathbf{K})}\mathrm{SKer}(f)
\right) =
\bigcap_{f\in\bigcup_{\mathbf{K}\in\mathcal{K}}\mathrm{Hom}(\mathbf{A},\mathbf{K})}\mathrm{Sch}(f).
$$
From now on, $\mathbf{T}_{\boldsymbol{\mathcal{K}}}(\mathbf{A})$ stands for $\bigcap_{f\in\bigcup_{\mathbf{K}\in\mathcal{K}}\mathrm{Hom}(\mathbf{A},\mathbf{K})}\mathbf{Sch}(f)/{\equiv^{\mathcal{K}}_{\mathbf{A}}}$. Then $\mathbf{T}_{\boldsymbol{\mathcal{K}}}(\mathbf{A})$ together with the restriction to $\mathbf{A}$ of the  projection from $\bigcap_{f\in\bigcup_{\mathbf{K}\in\mathcal{K}}\mathrm{Hom}(\mathbf{A},\mathbf{K})}\mathbf{Sch}(f)$ to $\mathbf{T}_{\boldsymbol{\mathcal{K}}}(\mathbf{A})$ is a universal morphism from $\mathbf{A}$ to $\mathrm{In}_{\boldsymbol{\mathcal{K}}}$. On the set $\bigcup_{\mathbf{K}\in\mathcal{K}}\mathrm{Hom}(\mathbf{A},\mathbf{K})$ let $\equiv$ be the equivalence relation defined, for every $f$, $g\in \bigcup_{\mathbf{K}\in\mathcal{K}}\mathrm{Hom}(\mathbf{A},\mathbf{K})$ as:
$$
f\equiv g\,\,\, \text{if and only if}\,\,\, \mathrm{SKer}(f) = \mathrm{SKer}(g),
$$
and let $\mathcal{H}$ be a subset of $\bigcup_{\mathbf{K}\in\mathcal{K}}\mathrm{Hom}(\mathbf{A},\mathbf{K})$ such that, for every $f\in \bigcup_{\mathbf{K}\in\mathcal{K}}\mathrm{Hom}(\mathbf{A},\mathbf{K})$, $\mathrm{card}(\mathcal{H}\cap [f]_{\equiv}) = 1$. Then $\mathbf{T}_{\boldsymbol{\mathcal{K}}}(\mathbf{A})$ is isomorphic to a subalgebra of $\prod_{f\in \mathcal{H}}\mathbf{Sch}(f)/ \mathrm{SKer}(f)$.

\textsf{Second proof (classical)}.\hfill

Let $\mathbf{A}$ be a partial $\Sigma$-algebra with $(\mathrm{card}(A_{s}))_{s\in S} = (\mathfrak{n}_{s})_{s\in S}$. Then, by Lemma~\ref{LCardinal}, there exists a cardinal $\mathfrak{m}_{\mathbf{A}}$ such that $(\mathrm{card}(B_{s}))_{s\in S}\leq (\mathfrak{m}_{\mathbf{A}})_{s\in S}$, for all partial $\Sigma$-algebras $\mathbf{B}$ generated by $S$-sorted sets $X$ such that $(\mathrm{card}(X_{s}))_{s\in S}\leq (\mathfrak{n}_{s})_{s\in S}$.

Let $\mathcal{K}_{\mathfrak{m}_{\mathbf{A}}}$ be the (small) set $\{\mathbf{K}\in \mathcal{K}\mid (\mathrm{card}(K_{s}))_{s\in S}\leq (\mathfrak{m}_{\mathbf{A}})_{s\in S}\}$. For $\mathbf{K}\in \mathcal{K}_{\mathfrak{m}_{\mathbf{A}}}$, let $(\mathbf{K}^{\mathrm{Hom}(\mathbf{A},\mathbf{K})},(\mathrm{pr}^{f})_{f\in \mathrm{Hom}(\mathbf{A},\mathbf{K})})$ be the product of $(\mathbf{K})_{f\in \mathrm{Hom}(\mathbf{A},\mathbf{K})}$. By the assumption on $\mathcal{K}$, $\mathbf{K}^{\mathrm{Hom}(\mathbf{A},\mathbf{K})}\in \mathcal{K}$. Then we denote by $\mathrm{p}^{\mathbf{K}}$ the unique homomorphism $\left<f\right>_{f\in\mathrm{Hom}(\mathbf{A},\mathbf{K})}$ from $\mathbf{A}$ to $\mathbf{K}^{\mathrm{Hom}(\mathbf{A},\mathbf{K})}$ such that, for every $f\in\mathrm{Hom}(\mathbf{A},\mathbf{K})$, $\mathrm{pr}^{f}\circ \left<f\right>_{f\in\mathrm{Hom}(\mathbf{A},\mathbf{K})} = f$.

Next, let $(\prod_{\mathbf{K}\in \mathcal{K}_{\mathfrak{m}_{\mathbf{A}}}}\mathbf{K}^{\mathrm{Hom}(\mathbf{A},\mathbf{K})},(\mathrm{pr}^{\mathbf{K}})_{\mathbf{K}\in \mathcal{K}_{\mathfrak{m}_{\mathbf{A}}}})$ be the product of $(\mathbf{K}^{\mathrm{Hom}(\mathbf{A},\mathbf{K})})_{\mathbf{K}\in \mathcal{K}_{\mathfrak{m}_{\mathbf{A}}}}$.  By the assumption on $\mathcal{K}$, $\prod_{\mathbf{K}\in \mathcal{K}_{\mathfrak{m}_{\mathbf{A}}}}\mathbf{K}^{\mathrm{Hom}(\mathbf{A},\mathbf{K})}\in \mathcal{K}$. Then we denote by $\mathrm{q}$ the unique homomorphism $\left<\mathrm{p}^{\mathbf{K}}\right>_{\mathbf{K}\in \mathcal{K}_{\mathfrak{m}_{\mathbf{A}}}}$ from $\mathbf{A}$ to
$\prod_{\mathbf{K}\in \mathcal{K}_{\mathfrak{m}_{\mathbf{A}}}}\mathbf{K}^{\mathrm{Hom}(\mathbf{A},\mathbf{K})}$ such that, for every $\mathbf{K}\in\mathcal{K}_{\mathfrak{m}_{\mathbf{A}}}$, $\mathrm{pr}^{\mathbf{K}}\circ \left<\mathrm{p}^{\mathbf{K}}\right>_{\mathbf{K}\in \mathcal{K}_{\mathfrak{m}_{\mathbf{A}}}} = \mathrm{p}^{\mathbf{K}}$. Let $\mathbf{T}_{\boldsymbol{\mathcal{K}}}(\mathbf{A})$ be the partial $\Sigma$-algebra induced by the subalgebra of $\prod_{\mathbf{K}\in \mathcal{K}_{\mathfrak{m}_{\mathbf{A}}}}\mathbf{K}^{\mathrm{Hom}(\mathbf{A},\mathbf{K})}$ generated by $\mathrm{q}[A]$, and $\theta^{\mathbf{A}}$ the unique epimorphism from $\mathbf{A}$ to $\mathbf{T}_{\boldsymbol{\mathcal{K}}}(\mathbf{A})$ such that $\mathrm{q} = \mathrm{in}^{\mathbf{T}_{\boldsymbol{\mathcal{K}}}(\mathbf{A})}\circ \theta^{\mathbf{A}}$. By the assumption on $\mathcal{K}$, $\mathbf{T}_{\boldsymbol{\mathcal{K}}}(\mathbf{A})\in \mathcal{K}$.

Let $\mathbf{B}$ be a partial $\Sigma$-algebra belonging to $\mathcal{K}$ and $g$ a homomorphism from $\mathbf{A}$ to $\mathbf{B}$. Let $\mathbf{Sg}_{\mathbf{B}}(g[A])$ be the partial $\Sigma$-algebra induced by the subalgebra $\mathrm{Sg}_{\mathbf{B}}(g[A])$ of $\mathbf{B}$ generated by $g[A]$, and $g^{\mathrm{e}}$ the unique epimorphism from $\mathbf{A}$ to $\mathbf{Sg}_{\mathbf{B}}(g[A])$ such that
$g = \mathrm{in}^{\mathbf{Sg}_{\mathbf{B}}(g[A])}\circ g^{\mathrm{e}}$.  By the assumption on $\mathcal{K}$, $\mathbf{Sg}_{\mathbf{B}}(g[A])\in \mathcal{K}$. Moreover, $(\mathrm{card}(g_{s}[A_{s}]))_{s\in S}\leq (\mathrm{card}(A_{s}))_{s\in S}$. Thus $(\mathrm{card}(\mathrm{Sg}_{\mathbf{B}}(g[A])_{s})_{s\in S}\leq (\mathfrak{m}_{\mathbf{A}})_{s\in S}$. Therefore $\mathbf{Sg}_{\mathbf{B}}(g[A])\in \mathcal{K}_{\mathfrak{m}_{\mathbf{A}}}$. Then
$$
\mathrm{in}^{\mathbf{Sg}_{\mathbf{B}}(g[A])}\circ\mathrm{pr}^{g^{\mathrm{e}}}\circ\mathrm{pr}^{\mathbf{Sg}_{\mathbf{B}}(g[A])}\circ\mathrm{in}^{\mathbf{T}_{\boldsymbol{\mathcal{K}}}(\mathbf{A})},
$$
denoted by $g^{\mathsf{p}}$, is the unique homomorphism from $\mathbf{T}_{\boldsymbol{\mathcal{K}}}(\mathbf{A})$ to $\mathbf{B}$ such that  $$g = g^{\mathsf{p}}\circ \theta^{\mathbf{A}}.$$

This fact is depicted in Figure~\ref{FFreeK}.
\end{proof}

\begin{figure}
$$\xymatrix@C=40pt@R=40pt{
{} & \mathbf{B} & {} \\
\mathbf{Sg}_{\mathbf{B}}(g[A])
\ar[ru]^-{\mathrm{in}^{\mathbf{Sg}_{\mathbf{B}}(g[A])}} &
\mathbf{A}
\ar[u]^-{g}
\ar[l]_-{g^{\mathrm{e}}}
\ar[r]^-{\theta^{\mathbf{A}}}
\ar[rd]^-{\mathrm{q}}
\ar[d]^-{\mathrm{p}^{\mathbf{Sg}_{\mathbf{B}}(g[A])}} &
\mathbf{T}_{\boldsymbol{\mathcal{K}}}(\mathbf{A})
\ar[ul]_-{g^{\mathsf{p}}}
\ar[d]^-{\mathrm{in}^{\mathbf{T}_{\boldsymbol{\mathcal{K}}}(\mathbf{A})}} \\
{} &
\mathbf{Sg}_{\mathbf{B}}(g[A])^{\mathrm{Hom}(\mathbf{A},\mathbf{Sg}_{\mathbf{B}}(g[A]))}
\ar[lu]^-{\mathrm{pr}^{g^{\mathrm{e}}}}
&
\prod_{\mathbf{K}\in \mathcal{K}_{\mathfrak{m}_{\mathbf{A}}}}\mathbf{K}^{\mathrm{Hom}(\mathbf{A},\mathbf{K})}
\ar[l]^-{\mathrm{pr}^{\mathbf{Sg}_{\mathbf{B}}(g[A])}}
   }
$$
\caption{The free $\boldsymbol{\mathcal{K}}$-algebra.}
\label{FFreeK}
\end{figure}

\begin{remark}
For $\theta^{\mathbf{A}}$ we have $\mathrm{Ker}(\theta^{\mathbf{A}}) = \bigcap_{\mathbf{K\in \mathcal{K}}}\left(
\bigcap_{f\in\mathrm{Hom}(\mathbf{A},\mathbf{K})}\mathrm{Ker}(f)\right)$.
\end{remark}

\begin{remark}
The second proof in the just stated theorem is interesting because it, essentially, is the model of the proof of the General Adjoint Functor of P. Freyd (see, e.g., \cite{em76} for the proof of Freyd's theorem). 

It seems to us appropriate to point out that in order to prove, via Freyd's theorem, that the inclusion functor $\mathrm{In}_{\boldsymbol{\mathcal{K}}}$ from $\boldsymbol{\mathcal{K}}$ to $\mathsf{PAlg}(\Sigma)$ has a left adjoint, it is sufficient to verify the following conditions: 
\begin{enumerate}
\item $\boldsymbol{\mathcal{K}}$ is $\boldsymbol{\mathcal{U}}$-locally small (i.e., has $\boldsymbol{\mathcal{U}}$-small hom sets),
\item $\boldsymbol{\mathcal{K}}$ has $\boldsymbol{\mathcal{U}}$-small limits (i.e., has $\boldsymbol{\mathcal{U}}$-small products and equalizers),
\item $\mathrm{In}_{\boldsymbol{\mathcal{K}}}$ preserves $\boldsymbol{\mathcal{U}}$-small products and equalizers and
\item $\mathrm{In}_{\boldsymbol{\mathcal{K}}}$ sastisfies the solution set condition at $\mathbf{K}$, for every $\mathbf{K}$ in $\mathcal{K}$.
\end{enumerate}
The last condition, which is the most convoluted of them all, is settled in the above mentioned second proof, the remaining conditions are obviously satisfied.
\end{remark}

\begin{corollary}
Let $\mathcal{K}$ be a set of partial $\Sigma$-algebras closed under the operators $\mathrm{I}$, $\mathrm{S}$, and $\mathrm{P}$. Then the functor $\mathbf{T}_{\Sigma,\boldsymbol{\mathcal{K}}} = \mathbf{T}_{\boldsymbol{\mathcal{K}}}\circ \mathbf{D}_{\Sigma}$ from $\mathsf{Set}^{S}$ to $\boldsymbol{\mathcal{K}}$ is a left adjoint of the functor $\mathrm{G}_{\Sigma}\circ \mathrm{In}_{\boldsymbol{\mathcal{K}}}$ from $\boldsymbol{\mathcal{K}}$ to $\mathsf{Set}^{S}$.
\end{corollary}

\begin{corollary}
If $\mathcal{K}$ is a finitary $\mathrm{E}$-variety, a finitary $\mathrm{ECE}$-variety, or a finitary $\mathrm{QE}$-variety, then the inclusion functor $\mathrm{In}_{\boldsymbol{\mathcal{K}}}$ from $\boldsymbol{\mathcal{K}}$ to $\mathsf{PAlg}(\Sigma)$ has a left adjoint, and the functor $\mathrm{G}_{\Sigma}\circ \mathrm{In}_{\boldsymbol{\mathcal{K}}}$ from $\boldsymbol{\mathcal{K}}$ to $\mathsf{Set}^{S}$ has a left adjoint.
\end{corollary}

Let us conclude this section by recalling the following remark by Burmeister in~\cite{pb86}, on p.~68, on the higher complexity of the generation of algebraic structures in the case of partial algebras as compared to that of  total algebras: 
\begin{quote}
Especially for partial algebras, a combination of these two principles: stepwise generation of the elements and the structure on one side and stepwise enlarging of the congruence relations, form an important tool in connection with the generation of algebraic structures. In the theory of total algebras one usually already starts with a total ``free'' algebra and one only has to generate an appropriate congruence, while for partial algebras one additionally has to construct the elements and the admissible structure.
\end{quote}
This remark of Burmeister, as the reader will have the opportunity to verify, will be widely corroborated throughout this paper.

\chapter{Adjoint functors, fibrations and the Grothendieck construction}

This chapter reviews the concepts of adjoint functor, equivalence between categories, split fibration, split indexed category, and the Grothendieck construction for the reader's convenience.
For a full treatment of the topics of this section we refer the reader to \cite{BW85, BW12, Gro71, bj99, w92}.

\section{Adjoint functors and equivalences}

We start recalling the notion of adjointness and state a characterization of it by means of the notion of  universal morphism from an object of a category $\mathsf{X}$ to a functor $G$ from a category $\mathsf{A}$ to $\mathsf{X}$. Before doing that we point out that we let ``$\ast$'' and ``$\circ$'' stand for the horizontal and vertical composition of natural transformations, respectively; that we denote by $\mathrm{Id}_{\mathsf{X}}$ the identity functor at $\mathsf{X}$; and that we denote by $\mathrm{id}_{G}$ the identity natural transformation at $G\colon \mathsf{A}\mor \mathsf{X}$.

\begin{definition}
An \emph{adjointness} is a $6$-tuple $(\mathsf{A},\mathsf{X},G,F,\eta,\varepsilon)$ where $\mathsf{A}$ and $\mathsf{X}$ are categories, $G$ a functor from $\mathsf{A}$ to $\mathsf{X}$, $F$ a functor from 
$\mathsf{X}$ to $\mathsf{A}$, $\eta$ a natural transformation from $\mathrm{Id}_{\mathsf{X}}$ to $G\circ F$ and $\varepsilon$ a natural transformation from $F\circ G$ to $\mathrm{Id}_{\mathsf{A}}$ subject to the so called triangular identities
$$
(G\ast \varepsilon)\circ (\eta\ast G) = \mathrm{id}_{G} \text{ and }
(\varepsilon\ast F)\circ (F\ast\eta) = \mathrm{id}_{F}.
$$
Thus, for every object $a$ of $\mathsf{A}$, the first equation asserts that 
the morphism $\eta_{G(a)}$ from $G(a)$ to $G(F(G(a)))$ when composed with the morphism $G(\varepsilon_{a})$ from $G(F(G(a)))$ to $G(a)$ is $\mathrm{id}_{G(a)}$, the identity morphism at $G(a)$.
If $(\mathsf{A},\mathsf{X},G,F,\eta,\varepsilon)$ is an adjointness, then we will say that $F$ is a \emph{left adjoint} of $G$, that $G$ is a \emph{right adjoint} of $F$, that $\eta$ is the \emph{unit} of the adjunction and that $\varepsilon$ is the \emph{counit} of the adjunction. Given two categories $\mathsf{A}$ and $\mathsf{X}$, it is also customary to speak of an \emph{adjunction from} $\mathsf{X}$ \emph{to} 
$\mathsf{A}$, meaning by it a $4$-tuple $(G,F,\eta,\varepsilon)$ such that 
$(\mathsf{A},\mathsf{X},G,F,\eta,\varepsilon)$ is an adjointness. In such a case, for brevity, and when there is no danger of confusion, we drop the dependence on $\eta$ and $\varepsilon$ and write just $F\dashv G$ for $(G,F,\eta,\varepsilon)$.
We will write $F\dashv G$ to indicate that $F$ is a left adjoint of $G$ or, equivalently, that $G$ is a right adjoint of $F$.
\end{definition}

The following proposition is very useful because the second part of it allows to construct an adjunction whenever we have a universal morphism from every object of a category to a functor. 

\begin{proposition}\label{ChAd}
Let $G$ be a functor from $\mathsf{A}$ to $\mathsf{X}$. Then the following two conditions on $G$ are equivalent:
\begin{enumerate}
\item There exists an adjointness $(\mathsf{A},\mathsf{X},G,F,\eta,\varepsilon)$.
\item For every object $x$ of $\mathsf{X}$ there exists an object $F(x)$ of $\mathsf{A}$ and a morphism $\eta_{x}$ from $x$ to $G(F(x))$ such that, for every object $b$ of $\mathsf{A}$ and every morphism $f$ from $x$ to $G(b)$, there exists a unique morphism $f^{\sharp}$ from $F(x)$ to $b$ such that 
$f = G(f^{\sharp})\circ \eta_{F(x)}$.
\end{enumerate}
\end{proposition}

\begin{definition}\label{MorAdj}
Let $(G,F,\eta,\varepsilon)$ be an adjunction from $\mathsf{X}$ to $\mathsf{A}$ and 
$(G',F',\eta',\varepsilon')$ an adjunction from $\mathsf{X}'$ to $\mathsf{A}'$. Then a \emph{map of adjunctions} from $(G,F,\eta,\varepsilon)$ to $(G',F',\eta',\varepsilon')$ is a pair $(K,L)$ of functors with $K\colon \mathsf{A}\mor \mathsf{A}'$, $L\colon \mathsf{X}\mor \mathsf{X}'$ such that 
$L\circ G = G'\circ K$, $K\circ F = F'\circ L$, $L\ast \eta = \eta'\ast L$ and $\varepsilon'\ast K = K\ast \varepsilon$. When no confusion will arise, we will speak of a \emph{map of adjunctions} from 
$F\dashv G$ to $F'\dashv G'$. 
If $\mathsf{A}$ is a subcategory of $\mathsf{A}'$, $\mathsf{X}$ a subcategory of 
$\mathsf{X}'$ and $K$, $L$ the respective canonical inclusions, then we will say that $(G,F,\eta,\varepsilon)$ is a \emph{subadjunction} of $(G',F',\eta',\varepsilon')$ or that $F\dashv G$ is a \emph{subadjunction} of $F'\dashv G'$. 
\end{definition}

\begin{definition}
Let $G$ be a functor from $\mathsf{A}$ to $\mathsf{X}$. We will say that $G$ is (1) \emph{faithful} if, for every $a$, $b\in \mathrm{Ob}(\mathsf{A})$ and every $f$, $g\in \mathrm{Hom}(a,b)$, if $G(f) = G(g)$, then $f = g$; (2) \emph{full} if, for every $a$, $b\in \mathrm{Ob}(\mathsf{A})$ and every 
$u\in \mathrm{Hom}(F(a),F(b))$ there exists an $f\in \mathrm{Hom}(a,b)$ such that $F(f) = u$; (3)  \emph{essentially surjective} if, for every $x\in \mathrm{Ob}(\mathsf{X})$ there exists an $a\in \mathrm{Ob}(\mathsf{A})$ such that $F(a)$ is isomorphic to $x$; and (4) an \emph{equivalence} if it is full, faithful and essentially surjective.
\end{definition}

Next, we set out the key facts about equivalence.

\begin{proposition}\label{Eqvt}
Let $G$ be a functor from $\mathsf{A}$ to $\mathsf{X}$. Then the following three  conditions on $G$ are equivalent:
\begin{enumerate}
\item $G$ is an equivalence of categories.
\item There exists a functor $F$ from $\mathsf{X}$ to $\mathsf{A}$ and natural transformations $\eta$ from $\mathrm{Id}_{\mathsf{X}}$ to $G\circ F$ and $\varepsilon$ from $F\circ G$ to $\mathrm{Id}_{\mathsf{A}}$ such that $(\mathsf{A},\mathsf{X},G,F,\eta,\varepsilon)$ is an adjunction from $\mathsf{X}$ to $\mathsf{A}$ and the unit $\eta$ and the counit $\varepsilon$ are natural isomorphisms.
\item There exists a functor $F$ from $\mathsf{X}$ to $\mathsf{A}$ and natural isomorphisms $F\circ G\cong \mathrm{Id}_{\mathsf{A}}$ and 
$G\circ F\cong \mathrm{Id}_{\mathsf{X}}$. 
\end{enumerate}
\end{proposition}

\begin{definition}\label{lali}
A functor $F$ from $\mathsf{X}$ to $\mathsf{A}$ is said to be a \emph{left adjoint left inverse} of a functor $G$ from $\mathsf{A}$ to $\mathsf{X}$ if $F\dashv G$ and the counit is the identity.
\end{definition}

In the following proposition we recall that an adjunction $(G,F,\eta,\varepsilon)$ from a category $\mathsf{X}$ to another $\mathsf{A}$ determines an equivalence between full subcategories of fixed points of $\mathsf{A}$ and $\mathsf{X}$ (for more details we refer the reader to~\cite{PT91}). 
\begin{proposition}\label{UOp}
Let $(G,F,\eta,\varepsilon)$ be an adjunction from a category $\mathsf{X}$ to another $\mathsf{A}$. Then the categories $\mathsf{Fix}(\eta)$ and $\mathsf{Fix}(\varepsilon)$, where $\mathsf{Fix}(\eta)$ is the full subcategory of $\mathsf{X}$ determined by 
$$
\{x\in \mathrm{Ob}(\mathsf{X})\mid \eta_{x}\text{ is an isomorphism}\}
$$ 
and 
$\mathsf{Fix}(\varepsilon)$ the full subcategory of $\mathsf{A}$ determined by
$$
\{a\in \mathrm{Ob}(\mathsf{A})\mid \varepsilon_{a}\text{ is an isomorphism}\},
$$ 
are equivalent, i.e., there exists a functor from $\mathsf{Fix}(\eta)$ to $\mathsf{Fix}(\varepsilon)$ which is an equivalence (one can say, following~\cite{Lk81}, that the equivalence between $\mathsf{Fix}(\eta)$ and 
$\mathsf{Fix}(\varepsilon)$ is the ``unity of opposites'').
\end{proposition}

%

\section{The Grothendieck construction}
In what follows we begin by defining the notions of cartesian morphism, fibration, morphism between fibrations, split fibration and morphism between split fibrations---and in doing so we follow the presentation by Barr and Wells in~\cite{BW12}, pp.~327--329, and by Jacobs in~\cite{bj99}, p.~73. After this we define the notions of split $\mathsf{I}$-indexed category, for a category $\mathsf{I}$, split indexed category and morphism between split indexed categories---also following the terminology used by Jacobs in~\cite{bj99}, p.~51. Then we define the Grothendieck construction (see~\cite{Gro71}) which is a functor from the category of split indexed categories to the category of split fibrations.

\begin{definition}
Let $\pi\colon \mathsf{E}\mor \mathsf{B}$ be a functor, $f\colon c\mor d$ a morphism of $\mathsf{B}$ and $y$ an object of $\mathsf{E}$ such that $\pi(y) = d$. We will say that a morphism $u\colon x\mor y$ of $\mathsf{E}$ is \emph{cartesian} for $f$ and $y$ if $\pi(u) = f$ and, for every object $z$ of $\mathsf{E}$, every morphism $v\colon z\mor y$ of $\mathsf{E}$ and every morphism $h\colon \pi(z)\mor c$ of $\mathsf{B}$ for which $f\circ h = \pi(v)$, there exists a unique morphism $w\colon z\mor x$ of $\mathsf{E}$ such that $u\circ w = v$ and $\pi(w) = h$. We let $\mathrm{Cart}(f,y)$ stand for the set of all cartesian morphisms for $f$ and $y$.

A functor $\pi\colon \mathsf{E}\mor \mathsf{B}$ is a \emph{fibration} if there exists a cartesian morphism for every morphism $f\colon c\mor d$ of $\mathsf{B}$ and every object $y$ of $\mathsf{E}$ for which $\pi(y) = d$. Sometimes a fibration $\pi\colon \mathsf{E}\mor \mathsf{B}$ is called a \emph{fibred category} or a \emph{category fibred over} $\mathsf{B}$ and denoted by $(\mathsf{E},\pi)$. A \emph{morphism of fibrations} from a fibration $\pi\colon \mathsf{E}\mor \mathsf{B}$ to another $\pi'\colon \mathsf{E}'\mor \mathsf{B}'$ is a pair of functors $H\colon \mathsf{E}\mor \mathsf{E}'$ and $K\colon \mathsf{B}\mor \mathsf{B}'$ such that $K\circ \pi = \pi'\circ H$ and $H$ sends cartesian morphisms in $\mathsf{E}$ to cartesian morphisms in $\mathsf{E}'$. We let $\mathsf{Fib}$ stand for the corresponding category.

Let $\mathrm{Mor}(\mathsf{B})\times_{\mathrm{Ob}(\mathsf{B})}\mathrm{Ob}(\mathsf{E})$ be the fibered product of $\mathrm{tg}_{\mathsf{B}}\colon \mathrm{Mor}(\mathsf{B})\mor \mathrm{Ob}(\mathsf{B})$, the mapping that sends a morphism of $\mathsf{B}$ to its target, and $\pi\colon \mathrm{Ob}(\mathsf{E})\mor \mathrm{Ob}(\mathsf{B})$, the object mapping of the functor $\pi$. Then a \emph{clivage} for a fibration $\pi\colon \mathsf{E}\mor \mathsf{B}$ is an element $\gamma$ of $\prod_{(f,y)\in \mathrm{Mor}(\mathsf{B})\times_{\mathrm{Ob}(\mathsf{B})}\mathrm{Ob}(\mathsf{E})}\mathrm{Cart}(f,y)$, i.e., a mapping that assigns to every pair $(f,y)$ in which $f\colon c\mor d$ is a morphism of $\mathsf{B}$ and $y$ an object of $\mathsf{E}$ for which $\pi(y) = d$ a morphism $\gamma_{f,y}$ of $\mathsf{E}$ that is cartesian for $f$ and $y$. The clivage  $\gamma$ is a \emph{splitting} of the fibration $\pi\colon \mathsf{E}\mor \mathsf{B}$ if it is such that (1) for every object $d$ of $\mathsf{B}$ and every object $y$ of $\mathsf{E}$ for which $\pi(y) d$, $\gamma_{\mathrm{id}_{d},y} = \mathrm{id}_{y}$ and (2)  for every pair of morphisms $f\colon c\mor d$ and $g\colon d\mor e$ of 
$\mathsf{B}$ and every pair of objects $y$ and $z$ of $\mathsf{E}$ if $\pi(z) = e$ and $y$ is the domain of $\gamma_{g,z}$, then $\gamma_{g,z}\circ\gamma_{f,y} = \gamma_{g\circ f,z}$.

A fibration $\pi\colon \mathsf{E}\mor \mathsf{B}$ is \emph{split} if it has a splitting. A \emph{morphism of split fibrations} from a split fibration $\pi\colon \mathsf{E}\mor \mathsf{B}$ to another $\pi'\colon \mathsf{E}'\mor \mathsf{B}'$ is a pair of functors $(H,K)$ as above where $H$ preserves the splitting up to equality (not up to isomorphism). We let $\mathsf{SFib}$ stand for the corresponding category.
\end{definition}

\begin{definition}
Let $\mathsf{I}$ be a category. A \emph{split} $\mathsf{I}$-\emph{indexed category} is a contravariant functor $F$ from $\mathsf{I}$ to $\mathsf{Cat}$. Given an object $i\in \mathrm{Ob}(\mathsf{I})$, we write $\mathsf{F}_{i}$ for the category $F(i)$, and given a morphism $\varphi\in \mathrm{Hom}_{\mathsf{I}}(i,j)$, we write $F_{\varphi}$ for the functor $F(\varphi)\colon \mathsf{F}_{j}\mor \mathsf{F}_{i}$. A \emph{split indexed category} is an ordered pair $(\mathsf{I},F)$ in which $\mathsf{I}$ is a category and $F$ a split $\mathsf{I}$-indexed category. A \emph{morphism of split indexed categories} from a split indexed category $(\mathsf{I},F)$ to another $(\mathsf{I}',F')$ is an ordered pair $(G,\eta)$, where $G$ is a functor from $\mathsf{I}$ to $\mathsf{I}'$ and $\eta$ a natural transformation from $F$ to $F'\circ G^{\mathrm{op}}$, where $G^{\mathrm{op}}$ is the dual of $G$. We let $\mathsf{SICat}$  stand for the corresponding category.
\end{definition}

\begin{proposition}
The category $\mathsf{SICat}$ is (small) complete.
\end{proposition}

We next define the object part of the Grothendieck construction which univocally  assigns to a split indexed category a split fibration---and in doing so we follow the notation used in~\cite{JY21} (Chapter~10) but adding the domain of $F$ as a superscript to indicate the contravariance of $F$.

\begin{definition}
Let $(\mathsf{I},F)$ be a split indexed category. Then the \emph{Grothen\-dieck construction} at $(\mathsf{I},F)$ is $(\int^{\mathsf{I}}F,\pi_{F})$, where $\int^{\mathsf{I}}F$ is the category defined as follows:
\begin{enumerate}
\item $\mathrm{Ob}(\int^{\mathsf{I}}F) = \bigcup_{i\in\mathrm{Ob}(\mathsf{I})}(\{i\}\times\mathrm{Ob}(\mathsf{F}_{i}))$.
\item For every $(i,x), (j,y)\in \mathrm{Ob}(\int^{\mathsf{I}}F)$, $\mathrm{Hom}_{\int^{\mathsf{I}}F}((i,x),(j,y))$ is the set of all ordered pairs $(\varphi,f)$, where $\varphi\in \mathrm{Hom}_{\mathsf{I}}(i,j)$ and $f\in \mathrm{Hom}_{\mathsf{F}_{i}}(x,F_{\varphi}(y))$.
\item For every $(i,x)\in \mathrm{Ob}(\int^{\mathsf{I}}F)$, the identity morphism at $(i,x)$ is given by $(\mathrm{id}_{i},\mathrm{id}_{x})$.
\item For every $(i,x), (j,y), (k,z)\in \mathrm{Ob}(\int^{\mathsf{I}}F)$, every $(\varphi,f)\colon (i,x)\mor (j,y)$, and every $(\psi,g)\colon (j,y)\mor (k,z)$, the composite morphism $(\psi,g)\circ (\varphi,f)$ from $(i,x)$ to $(k,z)$ is
     $$
      (\psi,g)\circ (\varphi,f) = (\psi\circ \varphi,F_{\psi}(g)\circ f)\colon (i,x)\mor (k,z).
     $$
Notice that $f\colon x\mor F_{\varphi}(y)$, $g\colon y\mor F_{\psi}(z)$, hence, taking into account that $F_{\varphi}$ is a functor from $\mathsf{F}_{j}$ to $\mathsf{F}_{i}$, $F_{\varphi}(g)\colon F_{\varphi}(y)\mor F_{\varphi}(F_{\psi}(z))$. Therefore $F_{\psi}(g)\circ f\colon x\mor F_{\varphi}(F_{\psi}(z))$.
\end{enumerate}
And $\pi_{F}$ the functor from $\int^{\mathsf{I}}F$ to $\mathsf{I}$ that sends an object $(i,x)$ of $\int^{\mathsf{I}}F$ to $i$ and a morphism $(\varphi,f)\colon (i,x)\mor (j,y)$ of $\int^{\mathsf{I}}F$ to $\varphi\colon i\mor j$. The functor $\pi_{F}$ is a split fibration and we call it the \emph{canonical split fibration} from $\int^{\mathsf{I}}F$ to $\mathsf{I}$ (determined by $(\mathsf{I},F)$).
\end{definition}

We next define the morphism part of the Grothendieck construction. But before doing so, we prove that every morphism between split indexed categories determines a morphism between their associated split fibrations.

\begin{proposition}
Let $(G,\eta)$ be a morphism from a split indexed category $(\mathsf{I},F)$ to another  $(\mathsf{I}',F')$. Then there exists a functor $\int^{G}\eta$ from $\int^{\mathsf{I}}F$ to $\int^{\mathsf{I}'}F'$ such that $G\circ \pi_{F} = \pi_{F'}\circ \int^{G}\eta$ 
, i.e., such that the following diagram:
$$\xymatrix{
\int^{\mathsf{I}}F\ar[d]_{\pi_{F}}
\ar[r]^-{\int^{G}\eta} & \int^{\mathsf{I}'}F'\ar[d]^{\pi_{F'}}  \\
\mathsf{I}\ar[r]_-{G} &
\mathsf{I}'
}
$$
commutes.
\end{proposition}

\begin{proof}
If $(i,x)$ is an object of $\int^{\mathsf{I}}F$, then, taking into account that $\eta_{i}$ is a functor from $\mathsf{F}_{i}$ to $\mathsf{F}'_{G(i)}$, we define $(\int^{G}\eta)(i,x)$ as $(G(i),\eta_{i}(x))$. On the other hand, taking into account that the following diagram:
$$\xymatrix{
\mathsf{F}_{i} \ar[r]^-{\eta_{i}} & \mathsf{F}'_{G(i)} \\
\mathsf{F}_{j} \ar[u]^-{F_{\varphi}} \ar[r]_-{\eta_{j}} &
\mathsf{F}'_{G(j)} \ar[u]_-{F'_{G(\varphi)}}
}
$$
commutes, if $(\varphi,f)$ is a morphism from $(i,x)$ to $(j,y)$, then we define the morphism $(\int^{G}\eta)(\varphi,f)$ from $(G(i),\eta_{i}(x))$ to $(G(j),\eta_{j}(y))$ as  $(G(\varphi),\eta_{i}(f))$ (notice that $f\colon x\mor F_{\varphi}(y)$ and $\eta_{i}(F_{\varphi}(y)) = F'_{G(\varphi)}(\eta_{j}(y))$).
\end{proof}

\begin{definition}
Let $(G,\eta)$ be a morphism from a split indexed category $(\mathsf{I},F)$ to another $(\mathsf{I}',F')$. Then the \emph{Grothendieck construction} at $(G,\eta)$ is the morphism $(\int^{G}\eta, G)$ from the split fibration $(\int^{\mathsf{I}}F,\pi_{F})$ to the split fibration $(\int^{\mathsf{I}'}F',\pi_{F'})$. 
\end{definition}

\begin{proposition}
The Grothendieck construction is as equicalence between $\mathsf{SICat}$ and  
$\mathsf{SFib}$
\end{proposition}

\begin{corollary}
The category $\mathsf{SFib}$ is (small) complete.
\end{corollary}


The following propositions are useful to prove the completeness 
and the cocompleteness of categories obtained by means of the Grothendieck construction.

\begin{proposition}\cite[Theorem~1, p.~247]{tbg91}\label{PCompl}
Let $\mathsf{I}$ be a category and $F$ a contravariant functor from $\mathsf{I}$ to $\mathsf{Cat}$. If 
$\mathsf{I}$ is complete, for every object $i$ of $\mathsf{I}$, $\mathsf{F}_{i}$ is complete and, for every objects $i$, $j$ of $\mathsf{I}$ and every morphism $\varphi\in \mathrm{Hom}_{\mathsf{I}}(i,j)$, the functor $F_{\varphi}$ from $\mathsf{F}_{j}$ to $\mathsf{F}_{i}$ is continuous, then $\int^{\mathsf{I}}F$ is complete.
\end{proposition}

\begin{proposition}\cite[Theorem~2, p.~250]{tbg91}\label{PCoCompl}
Let $\mathsf{I}$ be a category and $F$ a contravariant functor from $\mathsf{I}$ to $\mathsf{Cat}$. If 
$\mathsf{I}$ is cocomplete, for every object $i$ of $\mathsf{I}$, $\mathsf{F}_{i}$ is cocomplete and, for every objects $i$, $j$ of $\mathsf{I}$ and every morphism $\varphi\in \mathrm{Hom}_{\mathsf{I}}(i,j)$, the functor $F_{\varphi}$ from $\mathsf{F}_{j}$ to $\mathsf{F}_{i}$ has a left adjoint, then $\int^{\mathsf{I}}F$ is cocomplete.
\end{proposition}

The following proposition provides a sufficient condition for a functor obtained by means of the Grothendieck construction to have a left adjoint. This is surely folklore. however, we give a proof because we have been unable to track down a complete proof of it.

\begin{proposition}\label{GCLadj}
Let $(G,\eta)$ be a morphism from the split indexed category $(\mathsf{I},F)$ to the split indexed category $(\mathsf{I}',F')$. If (1) $G$ has a left adjoint $T$, (2), for every $i'\in\mathrm{Ob}(\mathsf{I}')$, $F'_{\alpha'_{i'}}\colon \mathsf{F}'_{G(T(i'))}\mor \mathsf{F}'_{i'}$, where $\alpha'_{i'}\colon i'\mor G(T(i'))$ is the value at $i'$ of the unit of the adjunction $T\dashv G$, has a left adjoint $T'_{\alpha'_{i'}}$, and (3), for every $i\in\mathrm{Ob}(\mathsf{I})$, $\eta_{i}\colon \mathsf{F}_{i}\mor \mathsf{F}'_{G(i)}$ has a left adjoint $\xi_{i}$, then the functor $\int^{G}\eta$ from $\int^{\mathsf{I}}F$ to $\int^{\mathsf{I}'}F'$ has a left adjoint.
\end{proposition}

\begin{proof}
Let $(i',x')$ be an object of $\int^{\mathsf{I}'}F'$. We want to show that there exists a universal morphim from $(i',x')$ to $\int^{G}\eta$, i.e., that there exists an object $(i,x)$ of $\int^{\mathsf{I}}F$ and a morphism $(\varphi',f')$ in $\int^{\mathsf{I}'}F'$ from $(i',x')$ to $(\int^{G}\eta)(i,x) = (G(i),\eta_{i}(x))$ (hence such that $\varphi'\colon i'\mor G(i)$ and $f'\colon x'\mor F'_{\varphi'}(\eta_{i}(x))$) such that, for every object $(j,y)$ of $\int^{\mathsf{I}}F$ and every morphism $(\psi',g')$ from $(i',x')$ to $(\int^{G}\eta)(j,y) = (G(j),\eta_{j}(y))$ in $\int^{\mathsf{I}'}F'$ (hence $\psi'\colon i'\mor G(j)$ and $g'\colon x'\mor F'_{\psi'}(\eta_{j}(y))$), there exists a unique morphism $(\psi,g)$ from $(i,x)$ to $(j,y)$ (hence $\psi\colon i\mor j$ and $g\colon x\mor F_{\psi}(y)$) in $\int^{I}F$, such that
$$
\textstyle
(\int^{G}\eta)(\psi,g)\circ (\varphi',f') = (G(\psi),\eta_{i}(g))\circ (\varphi',f') = (\psi',g').
$$
Let us notice that $G(\psi)\colon G(i)\mor G(j)$ and $\eta_{i}(g)\colon \eta_{i}(x)\mor F'_{G(\psi)}(\eta_{j}(y))$. However, since $F'_{G(\psi)}(\eta_{j}(y)) = \eta_{i}(F_{\psi}(y))$, we have that $\eta_{i}(g)\colon \eta_{i}(x)\mor \eta_{i}(F_{\psi}(y))$. Therefore, the first coordinate of $(G(\psi),\eta_{i}(g))\circ (\varphi',f')$ is $G(\psi)\circ \varphi'\colon i'\mor G(j)$,
while its second coordinate is $F'_{\varphi'}(\eta_{i}(g))\circ f'\colon x'\mor F'_{\varphi'}(F'_{G(\psi)}(\eta_{j}(y)))$. However, since $F'_{G(\psi)}(\eta_{j}(y)) = \eta_{i}(F_{\psi}(y))$, we have that $F'_{\varphi'}(\eta_{i}(g))\circ f'\colon x'\mor F'_{\varphi'}(\eta_{i}(F_{\psi}(y)))$.

One of our problems will be: to define $\psi$ and $\varphi'$ in such a way that $\psi' = G(\psi)\circ \varphi'$ (in order to obtain that $F'_{\varphi'}\circ F'_{G(\psi)} = F'_{G(\psi)\circ\varphi'} = F'_{\psi'}$).

Let $i$ be $T(i')$, $\varphi' = \alpha'_{i'}\colon i'\mor G(T(i'))$ (the value at $i'$ of the unit of the adjunction $T\dashv G$), $x = \xi_{T(i')}(T'_{\alpha'_{i'}}(x'))$, and $f' = \beta'_{x'}\colon x'\mor (F'_{\alpha'_{i'}}\circ\eta_{T(i')})(x)$ (the value at $x'$ of the unit of the adjunction $\xi_{T(i')}\circ T'_{\alpha'_{i'}}\dashv F'_{\alpha'_{i'}}\circ \eta_{T(i')}$).

Let $(\psi',g')$ be a morphism from $(i',x')$ to $(\int^{G}\eta)(j,y) = (G(j),\eta_{j}(y))$ in $\int^{\mathsf{I}'}F'$. Then, since $\psi'\colon i'\mor G(j)$ and $T\dashv G$, we take as morphism $\psi\colon i\mor j$ precisely $\psi'^{\sharp}\colon T(i')\mor j$, which is a morphism from $i$ to $j$ since, by definition, $i = T(i')$ (recall that $\mathrm{Hom}_{\mathsf{I}'}(i',G(j))$ is naturally isomorphic to $\mathrm{Hom}_{\mathsf{I}}(T(i'),j) = \mathrm{Hom}_{\mathsf{I}}(i,j)$, hence $\psi'^{\sharp}$ is the value at $\psi'$ of such a natural isomorphism). See below for the definition of $g\colon x\mor F_{\psi}(y)$ from $g'$ precisely as $g'^{\sharp}\colon x\mor F_{\psi}(y)$ (the unique morphism from $x$ to $F_{\psi}(y)$ such that $g' = F'_{\varphi'}(\eta_{i}(g'^{\sharp}))\circ f'$).

Then, since $G(i) = G(T(i'))$ and, by definition, $\psi = \psi'^{\sharp}$ and $\varphi'= \alpha'_{i'}$, we have that $\psi' = G(\psi)\circ \varphi' (= G(\psi'^{\sharp})\circ \alpha'_{i'})$. Hence the following diagram
$$\xymatrix{
\mathsf{F}'_{G(i)} \ar[r]^-{F'_{\varphi'}} & \mathsf{F}'_{i'}  \\
\mathsf{F}'_{G(j)} \ar[u]^-{F'_{G(\psi)}} \ar[ru]_-{F'_{\psi'}} & {}
}
$$
commutes. Therefore we can assert that  $F'_{\varphi'}(\eta_{i}(g))\circ f'\colon x'\mor F'_{\psi'}(\eta_{j}(y))$.

On the other hand, taking into account that $f' = \beta'_{x'}\colon x'\mor (F'_{\alpha'_{i'}}\circ\eta_{T(i')})(x)$ (the value at $x'$ of the unit of the adjunction $\xi_{T(i')}\circ T'_{\alpha'_{i'}}\dashv F'_{\alpha'_{i'}}\circ \eta_{T(i')}$), that $\varphi'$ denotes $\alpha'_{i'}$, that $g'\colon x'\mor F'_{\psi'}(\eta_{j}(y))$, that $F'_{\psi'}(\eta_{j}(y)) = F'_{\varphi'}(F'_{G(\psi)}(\eta_{j}(y)))$ (since the diagram above commutes, i.e., $F'_{\psi'} = F'_{\varphi'}\circ F'_{G(\psi)}$), and that $\eta_{i}\circ F_{\psi} = F'_{G(\psi)}\circ\eta_{j}$ (hence $F'_{\psi'}(\eta_{j}(y)) = F'_{\varphi'}(F'_{G(\psi)}(\eta_{j}(y))) = F'_{\varphi'}(\eta_{i}(F_{\psi}(y)))$), we have that there exists a unique morphism $g'^{\sharp}\colon x\mor F_{\psi}(y)$, denoted by $g$, such that the following diagram
$$\xymatrix{
x' \ar[r]^-{f'} \ar[rd]_-{g'} & F'_{\varphi'}(\eta_{i}(x)) \ar[d]^-{F'_{\varphi'}(\eta_{i}(g))} \\
{} &  F'_{\varphi'}(\eta_{i}(F_{\psi}(y)))
}
$$
commutes.

$$
\xymatrix@=5pc{
\mathsf{F}_{T(i')}\ar@<1.5ex>[r]^{\eta_{T(i')}}
\ar@{}[r]|{\uadj} &
\mathsf{F}'_{G(T(i'))}\ar@<1.5ex>[l]^{\xi_{T(i')}}\ar@<1.5ex>[r]^{F'_{\alpha'_{i'}}}
\ar@{}[r]|{\uadj} & \mathsf{F}'_{i'}
\ar@<1.5ex>[l]^{T'_{\alpha'_{i'}}}
}
$$
This completes the proof.
\end{proof}

\chapter{Higher-order categories}\label{S0D}

In this paper, for technical reasons, in addition to the usual many-sorted notions of $n$-category and of $\omega$-category, we will need to use single-sorted versions of them. In what follows, therefore, we will proceed to define both the many-sorted and the single-sorted presentation of higher-order categories and to prove their equivalence.

\section{Single-sorted presentation of higher-order categories}

\begin{definition}\label{DCat}
A \emph{single-sorted category} is an ordered pair $\mathsf{C} = (C,\xi)$ where $C$ is a set and $\xi = (\#,\mathrm{sc},\mathrm{tg})$ a structure of category on $C$, i.e., $\#\colon C\times C\dmor C$ is a partial operation of \emph{composition}, and $\mathrm{sc}$, $\mathrm{tg}$ are unary operations on $C$ which assign, to every morphism in $C$, its \emph{source} and its \emph{target}, respectively, subject to satisfy the following conditions:
\begin{enumerate}
\item[(C1)]\label{C1} For every $f\in C$, 
\allowdisplaybreaks
\begin{align*}
\mathrm{sc}(\mathrm{sc}(f)) &= \mathrm{sc}(f),&
 \mathrm{sc}(\mathrm{tg}(f)) &= \mathrm{tg}(f),\\
\mathrm{tg}(\mathrm{sc}(f)) &= \mathrm{sc}(f),&
\mathrm{tg}(\mathrm{tg}(f)) &= \mathrm{tg}(f).
\end{align*}
\item[(C2)]\label{C2} For every $f,g\in C$, 
\allowdisplaybreaks
\begin{align*}
g\#f &\mbox{ is defined}&
 \mbox{if and only if}&&
  \mathrm{sc}(g)& = \mathrm{tg}(f).
\end{align*} 
\item[(C3)]\label{C3} For every $f,g\in C$, if $\mathrm{sc}(g) = \mathrm{tg}(f)$, then 
\allowdisplaybreaks
\begin{align*}
\mathrm{sc}(g\# f) &= \mathrm{sc}(f),
&
\mathrm{tg}(g\# f) &= \mathrm{tg}(g).
\end{align*}
\item[(C4)]\label{C4} For every $f\in C$, 
\allowdisplaybreaks
\begin{align*}
f\# \mathrm{sc}(f) &= f,
&
\mathrm{tg}(f)\# f &= f.
\end{align*}
\item[(C5)]\label{C5} For every $f,g,h\in C$, if $\mathrm{sc}(h)=\mathrm{tg}(g)$ and $\mathrm{sc}(g)=\mathrm{tg}(f)$, then 
\allowdisplaybreaks
\begin{align*}
h\# (g\# f) = (h\# g)\# f.
\end{align*}
\end{enumerate}

We will call the $f\in C$ such that $\mathrm{sc}(f) = f$ or, what is equivalent, $\mathrm{tg}(f) = f$,  \emph{identities} of $C$. 

Let $\mathsf{C}=(C,\xi)$ and $\mathsf{C}'=(C',\xi')$ be two single-sorted categories. A \emph{single-sorted functor} from $\mathsf{C}$ to $\mathsf{C}'$ is an ordered triple $(\mathsf{C},F,\mathsf{C}')$, denoted by
$F\colon \mathsf{C}\mor\mathsf{C}'$, where $F$ is a mapping from $C$ to $C'$ such that
\begin{itemize}
\item[(Ci)]\label{Ci} If $f\in C$, then
\allowdisplaybreaks
\begin{align*}
F(\mathrm{sc}(f))&=\mathrm{sc}'(F(f)),
&
F(\mathrm{tg}(f))&=\mathrm{tg}'(F(f)).
\end{align*}
\item[(Cii)]\label{Cii} If $f,g\in C$ are such that $\mathrm{sc}(g)=\mathrm{tg}(f)$, then
\allowdisplaybreaks
\begin{align*}
F(g\#f)=F(g)\#'F(f).
\end{align*}
\end{itemize}

We will say that a single-sorted category $\mathsf{C}$ is \emph{discrete} if every $f\in C$ is an identity (i.e., if it contains no non-identity morphisms).


Single-sorted functors can be composed. In fact, if $F$ is a single-sorted functor from $\mathsf{C}$ to $\mathsf{D}$ and $G$ is a single-sorted functor from $\mathsf{D}$ to $\mathsf{E}$, then $G\circ F$ is a single-sorted functor from $\mathsf{C}$ to $\mathsf{E}$, which we call the \emph{composite} of $F$ and $G$. This composition is associative. Moreover, for every single-sorted category $\mathsf{C}$, there is an \emph{identity} single-sorted functor $(\mathsf{C},\mathrm{Id}_{C},\mathsf{C})$, denoted by $\mathrm{Id}_{\mathsf{C}}\colon \mathsf{C}\mor\mathsf{C}$, where $\mathrm{Id}_{C}$ is the identity mapping  at $C$, which is the neutral element for the composition. We denote by $\mathsf{Cat}^{}$ the category consisting of single-sorted categories and single-sorted functors.
\end{definition}

Concerning the notion of single-sorted $2$-category, let us recall what Mac Lane in~\cite{sM98}, p.~280 said:  ``Similarly a $2$-category can be considered to be a single set $C$ considered as the set of $2$-cells. Then the previous $1$-cells (the arrows) and the $0$-cells (the objects) are just regarded as special ``degenerate'' $2$-cells. On the set $C$ of $2$-cells there are then two category structures, the ``horizontal'' structure $(\circ_{h},\mathrm{d}_{h},\mathrm{cd}_{h})$ and the ``vertical'' structure $(\circ_{v},\mathrm{d}_{v},\mathrm{cd}_{v})$. Each satisfies the axioms above for a category structure and in addition (1) every identity for the $0$-structure is an identity for the $1$-structure; (2) the two category structures commute with each other.''

We next define the notion of single-sorted $2$-category and in doing so we follow Johnson (see~\cite{J87}, p.~6).

\begin{definition}\label{D2Cat}
A \emph{single-sorted $2$-category} is an ordered pair $\mathsf{C} = (C,(\xi_{k})_{k\in 2})$
where $C$ is a set and, for every $k\in 2$, $\xi_{k}= (\#^{k},\mathrm{sc}^{k},\mathrm{tg}^{k})$ a structure of single-sorted category on $C$, such that the following conditions are satisfied:
\begin{enumerate}
\item[(2C1)]\label{2C1} For every $f\in C$,
\allowdisplaybreaks
\begin{align*}
\mathrm{sc}^{1}(\mathrm{sc}^{0}(f))&=\mathrm{sc}^{0}(f)
=\mathrm{sc}^{0}(\mathrm{sc}^{1}(f))=\mathrm{sc}^{0}(\mathrm{tg}^{1}(f)),\\
\mathrm{tg}^{1}(\mathrm{tg}^{0}(f))&=\mathrm{tg}^{0}(f)
=\mathrm{tg}^{0}(\mathrm{tg}^{1}(f))=\mathrm{tg}^{0}(\mathrm{sc}^{1}(f)).
\end{align*}
\item[(2C2)]\label{2C2} For every $f,g\in C$, if $\mathrm{sc}^{0}(g)=\mathrm{tg}^{0}(f)$, then
\allowdisplaybreaks
\begin{align*}
\mathrm{sc}^{1}(g{\#^{0}}f)&=\mathrm{sc}^{1}(g){\#^{0}}\mathrm{sc}^{1}(f),
&
\mathrm{tg}^{1}(g{\#^{0}}f)&=\mathrm{tg}^{1}(g){\#^{0}}\mathrm{tg}^{1}(f).
\end{align*}
\item[(2C3)]\label{2C3} For every $f,f',g,g'\in C$, if 
\allowdisplaybreaks
\begin{align*}
\mathrm{sc}^{1}(g')&=\mathrm{tg}^{1}(g);
&
\mathrm{sc}^{1}(f')&=\mathrm{tg}^{1}(f);
\\
\mathrm{sc}^{0}(g')&=\mathrm{tg}^{0}(f');
&
\mathrm{sc}^{0}(g)&=\mathrm{tg}^{0}(f),
\end{align*}
then
\allowdisplaybreaks
\begin{align*}
(g'\#^{0}f')\#^{1}(g\#^{0} f)
=
(g'\#^{1}g)\#^{0}(f'\#^{1}f)
.
\end{align*}
\end{enumerate}


We will use
\begin{center}
\begin{tikzpicture}
[ACliment/.style={-{To [angle'=45, length=5.75pt, width=4pt, round]}}
, scale=0.8, 
AClimentD/.style={double equal sign distance,
-implies
}
]
\node[] (u) at (0,0) [] {$\mathrm{sc}^{0}(f)$};
\node[] (v) at (4,0) [] {$\mathrm{tg}^{0}(f)$};
\draw[ACliment, bend left] (u) to node [above]	{$\mathrm{sc}^{1}(f)$} (v);
\draw[ACliment, bend right] (u) to node [below] {$\mathrm{tg}^{1}(f)$} (v);
\node[] (a) at (2,.5) [] {};
\node[] (b) at (2,-.5) [] {};
\draw[AClimentD]  (a) to node [right]{$f$} (b);
\end{tikzpicture}
\end{center}
to denote that $f\in C$.

We identify the family $(\xi_{k})_{k\in 2}$ with the ordered pair $(\xi_{0},\xi_{1})$.

Let $\mathsf{C}=(C,(\xi_{k})_{k\in 2})$ and $\mathsf{C}'=(C',(\xi'_{k})_{k\in 2})$ be two single-sorted $2$-categories. A \emph{single-sorted $2$-functor} from $\mathsf{C}$ to $\mathsf{C}'$ is an ordered triple $(\mathsf{C},F,\mathsf{C}')$, denoted by $F\colon \mathsf{C}\mor\mathsf{C}'$, where $F$ is a mapping from $C$ to $C'$ such that, for every $k\in 2$, $F$ is a single-sorted functor from $(C,\xi_{k})$ to $(C',\xi'_{k})$. 

Single-sorted $2$-functors can be composed. In fact, if $F$ is a single-sorted $2$-functor from $\mathsf{C}$ to $\mathsf{D}$ and $G$ is a single-sorted $2$-functor from $\mathsf{D}$ to $\mathsf{E}$, then $G\circ F$ is a single-sorted $2$-functor from $\mathsf{C}$ to $\mathsf{E}$, which we call the \emph{composite} of $F$ and $G$. This composition is associative. Moreover, for every single-sorted $2$-category $\mathsf{C}$, there is an \emph{identity} single-sorted $2$-functor $(\mathsf{C},\mathrm{Id}_{C},\mathsf{C})$, denoted by $\mathrm{Id}_{\mathsf{C}}\colon \mathsf{C}\mor\mathsf{C}$, where $\mathrm{Id}_{C}$ is the identity mapping  at $C$, which is the neutral element for the composition. We denote by $\mathsf{2Cat}^{}$ the category consisting of single-sorted 
$2$-categories and single-sorted $2$-functors.
\end{definition}

In the definition of \(\omega\)-category, we follow the approach taken by Brown and Higgins in \cite{BH81} and that of Al-Agl, Brown and Steiner in \cite{ABS00}.

\begin{definition}\label{DnCatSS} 
A \emph{single-sorted $n$-category}, for $n\in\mathbb{N}$ with $n\geq 2$, is an ordered pair $\mathsf{C} = (C,(\xi_{k})_{k\in n})$ where $C$ is a set and, for every $k\in n$, $\xi_{k}$ a structure of  single-sorted category on $C$ such that, for every $j,k\in n$ with $j<k$, the ordered pair $(C,(\xi_{i})_{i\in\{j,k\}})$ is a single-sorted $2$-category. For $j<k$, we identify $(\xi_{i})_{i\in \{j,k\}}$ with $(\xi_{j},\xi_{k})$. 

For \(k\in n\) we will say that \(f\in C\) is a \emph{\(k\)-cell} if \(f=\mathrm{sc}^{k}(f)\) or, equivalently, by item~(C1) in Definition~\ref{DCat}, if \(f=\mathrm{tg}^{k}(f)\). Then we denote by \(C_{k}\) the set of \(k\)-cells of \(\mathsf{C}\), i.e.,
\[
C_{k}=\{f\in C\mid f=\mathrm{sc}^{k}(f)\}.
\]
Moreover, we let \(C_{n}\) stand for \(C\) and, by analogy to the above, the elements of \(C_{n}\) will be called \emph{\(n\)-cells}.

A \emph{single-sorted $n$-functor} $F$ from $\mathsf{C}$ to $\mathsf{C}'$, denoted by
$F\colon \mathsf{C}\mor\mathsf{C}'$, is a mapping from $C$ to $C'$ such that, for every $k\in n$, $F$ is a single-sorted functor from $(C, \xi_{k})$ to $(C', \xi'_{k})$. 

Single-sorted $n$-functors can be composed. In fact, if $F$ is a single-sorted $n$-functor from $\mathsf{C}$ to $\mathsf{D}$ and $G$ is a single-sorted $n$-functor from $\mathsf{D}$ to $\mathsf{E}$, then $G\circ F$ is a single-sorted $n$-functor from $\mathsf{C}$ to $\mathsf{E}$, which we call the \emph{composite} of $F$ and $G$. This composition is associative. Moreover, for every single-sorted $n$-category $\mathsf{C}$, there is an \emph{identity} single-sorted $n$-functor $(\mathsf{C},\mathrm{Id}_{C},\mathsf{C})$, denoted by $\mathrm{Id}_{\mathsf{C}}\colon \mathsf{C}\mor\mathsf{C}$, where $\mathrm{Id}_{C}$ is the identity mapping  at $C$, which is the neutral element for the composition. We denote by $\mathsf{nCat}^{}$ the category consisting of single-sorted $n$-categories and single-sorted $n$-functors.

A \emph{single-sorted $\omega$-category} is an ordered pair $\mathsf{C} = (C,(\xi_{k})_{k\in \omega})$ where $C$ is a set and, for every $k\in \omega$, $\xi_{k}$ a structure of single-sorted category on $C$ such that:
\begin{enumerate}
\item[(\(\omega\)C1)]\label{wC1}
For every $j,k\in \omega$ with $j<k$, the ordered pair $(C,(\xi_{i})_{i\in\{j,k\}})$ is a single-sorted $2$-category. As above, we identify $(\xi_{i})_{i\in \{j,k\}}$ with $(\xi_{j},\xi_{k})$.
\item[(\(\omega\)C2)]\label{wC2}
For every \(f \in C\), there exists \(k\in\omega\) such that \(f=\mathrm{sc}^{k}(f)\). 
\end{enumerate}

Similarly to the \(n\)-categorial case, for \(k\in\omega\) we will say that \(f\in C\) is a \emph{\(k\)-cell} if \(f=\mathrm{sc}^{k}(f)\) or, equivalently, by item~(C1) in Definition~\ref{DCat}, if \(f=\mathrm{tg}^{k}(f)\). Then we denote by \(C_{k}\) the set of \(k\)-cells of \(\mathsf{C}\), i.e.,
\(
C_{k}=\{f\in C\mid f=\mathrm{sc}^{k}(f)\}.
\)
Let us note that condition~(\(\omega\)C2) is equivalent to the condition \(C=\bigcup_{k\in\omega}C_{k}\).

A \emph{single-sorted $\omega$-functor} $F$ from $\mathsf{C}$ to $\mathsf{C}'$, denoted by $F\colon \mathsf{C}\mor\mathsf{C}'$, is a mapping from $C$ to $C'$ such that, for every $k\in \omega$, $F$ is a single-sorted functor from $(C,\xi_{k})$ to $(C,\xi'_{k})$.

Single-sorted $\omega$-functors can be composed. In fact, if $F$ is a single-sorted $\omega$-functor from $\mathsf{C}$ to $\mathsf{D}$ and $G$ is a single-sorted $\omega$-functor from $\mathsf{D}$ to $\mathsf{E}$, then $G\circ F$ is a single-sorted $\omega$-functor from $\mathsf{C}$ to $\mathsf{E}$, which we call the \emph{composite} of $F$ and $G$. This composition is associative. Moreover, for every single-sorted $\omega$-category $\mathsf{C}$, there is an \emph{identity} single-sorted $\omega$-functor $(\mathsf{C},\mathrm{Id}_{C},\mathsf{C})$, denoted by $\mathrm{Id}_{\mathsf{C}}\colon \mathsf{C}\mor\mathsf{C}$, where $\mathrm{Id}_{C}$ is the identity mapping  at $C$, which is the neutral element for the composition. We denote by $\w\mathsf{Cat}$ the category consisting of single-sorted $\omega$-categories and single-sorted $\omega$-functors.
\end{definition}

\begin{remark}
 Note that, for $n\in\mathbb{N}$, a single-sorted $n$-category is a single-sorted $\omega$-category $\mathsf{C} = (C,(\#^{k},\mathrm{sc}^{k},\mathrm{tg}^{k})_{k\in \omega})$ such that for every $k\geq n$, the single-sorted category $(C,(\#^{k},\mathrm{sc}^{k},\mathrm{tg}^{k}))$ is discrete.   
\end{remark}

\section{
\texorpdfstring
{\(\w\mathsf{Cat}\) as the projective limit of \(((\bigcdot)\mathsf{Cat}, U)\)}
{Omega-Cat as a projective limit}
}

In this section, we define, for \(m,n\in\omega\) with \(n\leq m\), the single-sorted underlying functors \(U^{(n,m)}\) and \(U^{(n,\omega)}\). Moreover, we prove that the category \(\w\mathsf{Cat}\) of single-sorted \(\omega\)-categories with single-sorted \(\omega\)-functors is the projective limit of the projective system \(((\bigcdot)\mathsf{Cat},U)\) where the components of the projective system will be introduced in Remark~\ref{RUmnProjSys}. 

We begin by defining the notions of the underlying \(m\)-category and the underlying functors.

\begin{definition}\label{DUnderlying}
Let \(\mathsf{C} = (C, (\xi_{k})_{k\in n})\) be a single-sorted \(n\)-category with 
\[
(\xi_{k})_{k\in n} = (\#^{k}, \mathrm{sc}^{k}, \mathrm{tg}^{k})_{k\in n}.
\]

For \(m<n\), concerning \(C_{m}\), the set of \(m\)-cells of \(\mathsf{C}\), we define the following mappings.

For every \(0 \leq k < m\), we have \(\mathrm{sc}^{k}[C_{m}] \subseteq C_{m}\). Indeed, let \(f\) be an \(m\)-cell in \(C\), i.e., \(f = \mathrm{sc}^{m}(f)\). Then, by item~(2C1) in Definition~\ref{D2Cat}, it happens that
\[
\mathrm{sc}^{m}(\mathrm{sc}^{k}(f)) = \mathrm{sc}^{k}(f).
\]
Therefore, \(\mathrm{sc}^{k}\mathnormal{\upharpoonright}_{C_{m}}\) corestricts to \(C_{m}\).

For every \(0 \leq k < m\), we have \(\mathrm{tg}^{k}[C_{m}] \subseteq C_m\). Indeed, let \(f\) be an \(m\)-cell in \(C_{m}\), i.e., \(f = \mathrm{tg}^{m}(f)\). Then, by item~(2C1) in Definition~\ref{D2Cat}, it happens that
\[
\mathrm{tg}^{m}(\mathrm{tg}^{k}(f)) = \mathrm{tg}^{k}(f).
\]
Therefore, \(\mathrm{tg}^{k}\mathnormal{\upharpoonright}_{C_{m}}\) corestricts to \(C_{m}\).

Moreover, for every \(0 \leq k < m\), if \(f,g \in C_{m}\) are such that \(\mathrm{sc}^{k}(g) = \mathrm{tg}^{k}(f)\), then we have that the \(k\)-composite \(g \#^{k} f\) belongs to \(C_{m}\). Indeed, by item~(2C2) in Definition~\ref{D2Cat} and taking into account that \(f\) and 	\(g\) are \(m\)-cells, it happens that
\[
\mathrm{sc}^{m}(g \#^{k} f) = \mathrm{sc}^{m} (g) \#^{k} \mathrm{tg}^{m}(f) = g \#^{k} f.
\]
Therefore, \(\#^{k}\mathnormal{\upharpoonright}_{C_{m} \times C_{m}}\) corestricts to \(C_{m}\).

We will call the ordered pair
\[
(C_m, (\xi^{<m}_{k})_{k \in m})
\]
with \((\xi^{<m}_{k})_{k\in m} = (\#^k\bigr|^{C_{m}}_{C_{m}\times C_{m}}, \mathrm{sc}^{k}\bigr|^{C_m}_{C_m}, \mathrm{tg}^{k}\bigr|^{C_m}_{C_m})_{k \in m}\) the \emph{underlying single-sorted \(m\)-category} of \(\mathsf{C}\) and we will denote it by \(\mathsf{C}^{<m}\). Thus, it follows that \(\mathsf{C}^{<m}\) is a single-sorted \(m\)-category.

Moreover, if \(F\) is a single-sorted \(n\)-functor from \(\mathsf{C} = (C,(\xi_{k})_{k\in n})\) to \(\mathsf{C}' = (C',(\xi'_{k})_{k\in n})\). We have that \(F[C_{m}]\subseteq C'_{m}\). Indeed, let \(f\) be an \(m\)-cell in \(C\), i.e., \(f = \mathrm{sc}^{m}(f)\). Then, by item~(Ci) in Definition~\ref{DCat} it happens that
\[
\mathrm{sc}'^{m}(F(f)) = F(\mathrm{sc}^{m}(f)) = F(f).
\]
Therefore, \(F\mathnormal{\upharpoonright}_{C_{m}}\) corestricts to \(C'_{m}\). We will call the mapping \(F\bigr|^{C'_{m}}_{C_{m}}\) the \emph{underlying single-sorted \(m\)-functor} of \(F\) and we will denote it by \(F^{<m}\). Thus, it follows that \(F^{<m}\) is a single-sorted \(m\)-functor from \(\mathsf{C}^{<m}\) to \(\mathsf{C}'^{<m}\).

A similar construction applies if \(\mathsf{C}=(C,(\xi_{k})_{k\in\omega})\) is a \(\omega\)-category and \(F\) a \(\omega\)-functor.
\end{definition}

\begin{definition}\label{DUmn}
Let \(m,n\in\omega\) with \(m<n\). we let \(U^{(m,n)}\) stand for the covariant functor from \(\mathsf{nCat}\) to \(\mathsf{mCat}\) defined as follows:
\begin{enumerate}
\item
its object mapping sends a single-sorted \(n\)-category \(\mathsf{C}\) to its underlying single-sorted \(m\)-category \(\mathsf{C}^{<m}\); and
\item
its morphism mapping sends a single-sorted \(n\)-functor \(F\) from \(\mathsf{C}\) to \(\mathsf{C}'\) to its underlying single-sorted \(m\)-functor \(F^{<m}\) from \(\mathsf{C}^{<m}\) to \(\mathsf{C}'^{<m}\).
\end{enumerate}

Following a similar construction, for \(n\in\omega\) we define the covariant functor \(U^{(n,\omega)}\) from \(\w\mathsf{Cat}\) to \(\mathsf{nCat}\).

Moreover, we let \(U^{(n,n)}\) stand for \(\mathrm{Id}_{\mathsf{nCat}}\) the identity single-sorted functor on \(\mathsf{nCat}\).
\end{definition}

\begin{remark}\label{RUmnProjSys}
We let \((\bigcdot)\mathsf{Cat}\) stand for the family of single-sorted categories \((\mathsf{nCat})_{n\in\omega}\) and we let \(U\) stand for the family of functors \((U^{(m,n)})_{(m,n)\in\leq}\). Following Definition~\ref{DUnderlying}, for every \(m<n<p\), every \(p\)-category \(\mathsf{C}\) and every \(p\)-functor \(F\), it follows that \((\mathsf{C}^{<n})^{<m}=\mathsf{C}^{<m}\) and \((F^{<n})^{<m}=F^{<m}\), that is, 
\[
U^{(m,n)}\circ U^{(n,p)}=U^{(m,p)}.
\]
Therefore, we obtain the projective system \(((\bigcdot)\mathsf{Cat},U)\) in \(\mathsf{Cat}\) and its projective limit \(\varprojlim ((\bigcdot)\mathsf{Cat},U)\).
\end{remark}

Our next goal is to prove that such limit is isomorphic to \((\w\mathsf{Cat},U^{(\bigcdot,\omega)})\) with \(U^{(\bigcdot,\omega)}\) being the family of functors \((U^{(n,\omega)})_{n\in\omega}\). That the composition of functors \(U^{(m,n)}\circ U^{(n,\omega)}\) is equal to the functor \(U^{(m,\omega)}\) follows from the just stated consideration. Thus, all that remains to be proven is that \((\w\mathsf{Cat},U^{(\bigcdot,\omega)})\) satisfies the universal property.

\begin{proposition}
\label{PUnivPropSS}
Let \(\mathsf{L}\) be a category and, for every \(n\in\omega\), let \(F_{n}\) be a functor from \(\mathsf{L}\) to \(\mathsf{nCat}\). If, for every \(m\leq n\), \(F_{m} = U^{(m,n)} \circ F_{n}\), then there exists a unique functor \(\langle F_{k}\rangle_{k\in\omega}\) from \(\mathsf{L}\) to the category \(\w\mathsf{Cat}\) such that, for every \(n\in\omega\), \(F_{n} = U^{(n,\omega)} \circ \langle F_{k}\rangle_{k\in\omega}\). The diagram in Figure~\ref{FUnivPropSS} depicts this situation.
\begin{figure}
\[
\xymatrix{
& \mathsf{L}\ar@{-->}[d]^{\exists!\langle F_{k}\rangle_{k\in\omega}}\ar@/_30pt/[ddl]_{F_{m}}\ar@/^30pt/[ddr]^{F_{n}} &
\\
& \w\mathsf{Cat}\ar[dr]^{U^{(n,\omega)}}\ar[dl]_{U^{(m,\omega)}} &
\\
\mathsf{mCat} && \mathsf{nCat}\ar[ll]^{U^{(m,n)}}
}
\]
\caption{Universal property of \(\w\mathsf{Cat}\).}
\label{FUnivPropSS}
\end{figure}
\end{proposition}

\begin{proof}
Before proceeding any further let us fix some notation. For every object \(A\) in \(\mathsf{L}\), we will denote its image under the functor \(F_{n}\) by \(\mathsf{A}_{n}=(A_{n},(\xi^{\mathsf{A}_{n}}_{k})_{k\in n})\) with \((\xi^{\mathsf{A}_{n}}_{k})_{k\in n}=(\#^{k\mathsf{A}_{n}},\mathrm{sc}^{k\mathsf{A}_{n}},\mathrm{tg}^{k\mathsf{A}_{n}})_{k\in n}\). Let us remark that \(\mathsf{A}_{n}\) is a single-sorted \(n\)-category. Moreover, for every morphism \(\varphi\) from \(A\) to \(B\) in \(\mathsf{L}\), we will denote its image under the functor \(F_{n}\) by \(\varphi_{n}\). Let us remark that \(\varphi_{n}\) is a single-sorted \(n\)-functor from \(\mathsf{A}_{n}\) to \(\mathsf{B}_{n}\), thus a mapping from \(A_{n}\) to \(B_{n}\).

In order to define the functor \(\langle F_{k}\rangle_{k\in\omega}\), we will prove the following two claims.

\begin{claim}\label{CUnivPropSSA}
For every object \(A\) in \(\mathsf{L}\) we can uniquely assign a single-sorted \(\omega\)-category \(\mathsf{A}_{\omega}\) such that, for every \(n\in\omega\),
\[
\mathsf{A}_{\omega}^{<n}
=
\mathsf{A}_{n}.
\]
\end{claim}

Let us point out that, for every \(m, n \in \omega\) with \(m\leq n\), since \(F_{m}=U^{(m,n)}\circ F_{n}\), it follows that \(\mathsf{A}_{m}=\mathsf{A}_{n}^{<m}\). In particular, the following equalities hold
\begin{itemize}
\item[(CC1)]\label{CC1}
\(A_{m}=(A_{n})_{m}\) and \(A_{m}\subseteq A_{n}\); and 
\item[(CC2)]\label{CC2}
for every \(k\in m\), \(\xi^{\mathsf{A}_{m}}_{k}=(\xi^{\mathsf{A}_{n}}_{k})^{<m}\), that is,
\begin{align*}
\#^{k\mathsf{A}_{m}} &= \#^{k\mathsf{A}_{n}}\bigr|^{\mathsf{A}_{m}}_{\mathsf{A}_{m}\times \mathsf{A}_{m}}
&
\mathrm{tg}^{k\mathsf{A}_{m}} &= \mathrm{tg}^{k\mathsf{A}_{n}}\bigr|^{\mathsf{A}_{m}}_{\mathsf{A}_{m}}
&
\mathrm{sc}^{k\mathsf{A}_{m}} &= \mathrm{sc}^{k\mathsf{A}_{n}}\bigr|^{\mathsf{A}_{m}}_{\mathsf{A}_{m}}.
\end{align*}
\end{itemize}

Let \(A\) be an object in \(\mathsf{L}\). We define the ordered pair \(\mathsf{A}_{\omega}=(A_{\omega}, (\xi^{\mathsf{A}}_{k})_{k\in\omega})\) as follows:
\begin{enumerate}
\item
\(A_{\omega}\) is the set \(\bigcup_{k\in\omega}A_{k}\); and
\item
\((\xi^{\mathsf{A}}_{k})_{k\in\omega}=(\#^{k\mathsf{A}},\mathrm{sc}^{k\mathsf{A}},\mathrm{tg}^{k\mathsf{A}})_{k\in\omega}\) where, for every \(k\in\omega\),
\begin{enumerate}
\item
\(\mathrm{sc}^{k\mathsf{A}}\) and \(\mathrm{tg}^{k\mathsf{A}}\) are unary operations on \(A_{\omega}\) defined, for every \(f \in A_{\omega}\), as \(\mathrm{sc}^{k\mathsf{A}_{n_{k,f}}}(f)\) and \(\mathrm{tg}^{k\mathsf{A}_{n_{k,f}}}(f)\), respectively, where 
\[
n_{k,f}=\min\{ i\in\mathbb{N} \mid k+1\leq i, f\in A_{i}\}.
\]
Let us note that, from item~(CC2) above, it follows that, for every \(n\in\mathbb{N}\) with \(n_{k,f}\leq n\), \(\mathrm{sc}^{k\mathsf{A}}(f)=\mathrm{sc}^{k\mathsf{A}_{n}}f\) and \(\mathrm{tg}^{k\mathsf{A}}(f)=\mathrm{tg}^{k\mathsf{A}_{n}}f\). The index \(n_{k,f}\) should be read as the minimum natural number for which a \(k\)-th structural operation can be performed to the element \(f\).
\item
\(\#^{k\mathsf{A}}\) is a partial binary operation on \(A_{\omega}\) defined, for every \(f,g\in A_{\omega}\) with \(\mathrm{tg}^{k\mathsf{A}}(g)=\mathrm{sc}^{k\mathsf{A}}(f)\), as \(g\#^{k\mathsf{A}_{n_{k,f,g}}}f\) where
\[
n_{k,f,g} = \min\{ i\in\mathbb{N} \mid k+1\leq i, f,g\in A_{i} \}.
\]
Let us note that, from item~(CC1) above such a minimum exists. Moreover, from item~(CC2), it follows that, for every \(n\in\mathbb{N}\) with \(n_{k,f,g}\leq n\), \(g\#^{k\mathsf{A}}f=g\#^{k\mathsf{A}_{n}}f\). Similarly to \(n_{k,f}\), the index \(n_{k,f,g}\) should be read as the minimum natural number for which a \(k\)-th structural operation can be performed to the elements \(f\) and \(g\).
\end{enumerate}
\end{enumerate}

We want to prove that \(\mathsf{A}_{\omega}\) is a \(\omega\)-category. Thus, we must check that is satisfies the conditions stated in Definition~\ref{DnCatSS}.

We begin by proving that, for \(k\in\omega\), \((A_{\omega},\xi^{\mathsf{A}}_{k})\) is a single-sorted category. To this end, we consider the different items stated in Definition~\ref{DCat}.

\textsf{(C1)}
For every \(f\in A_{\omega}\), we define
\begin{enumerate}
\item
\(n_{k,f}=\min\{i\in\mathbb{N}\mid k+1\leq i, f\in A_{i}\}\), and
\item
\(n_{k,\mathrm{sc}(f)}=\min\{i\in\mathbb{N}\mid k+1\leq i, \mathrm{sc}^{k\mathsf{A}_{n_{k,f}}}(f)\in A_{i}\}\).
\end{enumerate}
Let us note that \(n_{k,\mathrm{sc}(f)}\leq n_{k,f}\).

The following chain of equalities holds
\begin{align*}
\mathrm{sc}^{k\mathsf{A}}(\mathrm{sc}^{k\mathsf{A}}(f))
&=
\mathrm{sc}^{k\mathsf{A}}(\mathrm{sc}^{k\mathsf{A}_{n_{k,f}}}(f))
\tag{1}
\\
&=
\mathrm{sc}^{k\mathsf{A}_{n_{k,\mathrm{sc}(f)}}}(\mathrm{sc}^{k\mathsf{A}_{n_{k,f}}}(f))
\tag{2}
\\
&=
\mathrm{sc}^{k\mathsf{A}_{n_{k,f}}}(\mathrm{sc}^{k\mathsf{A}_{n_{k,f}}}(f))
\tag{3}
\\
&=
\mathrm{sc}^{k\mathsf{A}_{n_{k,f}}}(f)
\tag{4}
\\
&=
\mathrm{sc}^{k\mathsf{A}}(f).
\tag{5}
\end{align*}

The first and second equality unravels the definition of \(\mathrm{sc}^{k\mathsf{A}}\); the third equality follows from the fact that \(n_{k,\mathrm{sc}(f)}\leq n_{k,f}\); the fourth equality follows from item~(C1) in Definition~\ref{DCat} since \((A_{n_{k,f}},\xi^{\mathsf{A}_{n_{k,f}}}_{k})\) is a single-sorted category; finally, the last equality recovers the definition of \(\mathrm{sc}^{k\mathsf{A}}\).

By a similar argument we also have that
\begin{align*}
\mathrm{sc}^{k\mathsf{A}}\circ\mathrm{tg}^{k\mathsf{A}}=\mathrm{tg}^{k\mathsf{A}};&
&\mathrm{tg}^{k\mathsf{A}}\circ\mathrm{sc}^{k\mathsf{A}}=\mathrm{sc}^{k\mathsf{A}};&
&\mathrm{tg}^{k\mathsf{A}}\circ\mathrm{tg}^{k\mathsf{A}}=\mathrm{tg}^{k\mathsf{A}}.
\end{align*}

\textsf{(C2)}
For every \(f,g\in A_{\omega}\), we define
\begin{enumerate}
\item
\(n_{k,f,g}=\min\{i\in\mathbb{N} \mid k+1\leq i, f,g\in A_{i}\}\),
\item
\(n_{k,f}=\min\{i\in\mathbb{N} \mid k+1\leq i, f\in A_{i}\}\), and
\item
\(n_{k,g}=\min\{i\in\mathbb{N} \mid k+1\leq i, g\in A_{i}\}\).
\end{enumerate}
Let us note that \(n_{k,f}\leq n_{k,f,g}\) and \(n_{k,g}\leq n_{k,f,g}\).

The following chain of equivalences holds
\allowdisplaybreaks
\begin{align*}
g\#^{k\mathsf{A}}f \mbox{ is defined}
&\Leftrightarrow
g\#^{k\mathsf{A}_{n_{k,f,g}}}f \mbox{ is defined}
\tag{1}
\\
&\Leftrightarrow
\mathrm{sc}^{k\mathsf{A}_{n_{k,f,g}}}(g)=\mathrm{tg}^{k\mathsf{A}_{n_{k,f,g}}}(f)
\tag{2}
\\
&\Leftrightarrow
\mathrm{sc}^{k\mathsf{A}_{n_{k,g}}}(g)=\mathrm{tg}^{k\mathsf{A}_{n_{k,f}}}(f)
\tag{3}
\\
&\Leftrightarrow
\mathrm{sc}^{k\mathsf{A}}(g)=\mathrm{tg}^{k\mathsf{A}}(f).
\tag{4}
\end{align*}

The first equivalence unravels the definition of \(\#^{k\mathsf{A}}\); the second equivalence follows from item~(C2) in Definition~\ref{DCat} since \((A_{n_{k,f,g}},\xi^{\mathsf{A}_{n_{k,f,g}}}_{k})\) is a single-sorted category; the third equivalence follows from the fact that \(n_{k,f}\leq n_{k,f,g}\) and \(n_{k,g}\leq n_{k,f,g}\); finally, the last equality follows from the definition of \(\mathrm{sc}^{k\mathsf{A}}\) and \(\mathrm{tg}^{k\mathsf{A}}\).

\textsf{(C3)}
For every \(f,g\in A_{\omega}\), if \(\mathrm{sc}^{k\mathsf{A}}(g)=\mathrm{tg}^{k\mathsf{A}}(f)\), then we define
\begin{enumerate}
\item
\(n_{k,f,g}=\min\{i\in\mathbb{N} \mid k+1\leq i, f,g\in A_{i}\}\),
\item
\(n_{k,g\#f}=\min\{i\in\mathbb{N} \mid k+1\leq i, g\#^{k\mathsf{A}_{n_{k,f,g}}}f\in A_{i}\}\), and
\item
\(n_{k,f}=\min\{i\in\mathbb{N} \mid k+1\leq i, f\in A_{i}\}\).
\end{enumerate}
Let us note that \(n_{k,g\#f}\leq n_{k,f,g}\) and \(n_{k,f}\leq n_{k,f,g}\).

The following chain of equalities holds
\begin{align*}
\mathrm{sc}^{k\mathsf{A}}(g\#^{k\mathsf{A}}f)
&=
\mathrm{sc}^{k\mathsf{A}}(g\#^{k\mathsf{A}_{n_{k,f,g}}}f)
\tag{1}
\\
&=
\mathrm{sc}^{k\mathsf{A}_{n_{k,g\#f}}}(g\#^{k\mathsf{A}_{n_{k,f,g}}}f)
\tag{2}
\\
&=
\mathrm{sc}^{k\mathsf{A}_{n_{k,f,g}}}(g\#^{k\mathsf{A}_{n_{k,f,g}}}f)
\tag{3}
\\
&=
\mathrm{sc}^{k\mathsf{A}_{n_{k,f,g}}}(f)
\tag{4}
\\
&=
\mathrm{sc}^{k\mathsf{A}_{n_{k,f}}}(f)
\tag{5}
\\
&=
\mathrm{sc}^{k\mathsf{A}}(f).
\tag{6}
\end{align*}

The first and second equalities unravel the definitions of \(\#^{k\mathsf{A}}\) and \(\mathrm{sc}^{k\mathsf{A}}\) respectively; the third equality follows from the fact that \(n_{k,g\#f}\leq n_{k,f,g}\); the fourth  equality follows from item~(C3) in Definition~\ref{DCat} since \((A_{n_{k,f,g}},\xi^{\mathsf{A}_{n_{k,f,g}}}_{k})\) is a single-sorted category; the fifth equality follows from the fact that \(n_{k,f}\leq n_{k,f,g}\); finally, the last equality recovers the definition of \(\mathrm{sc}^{k\mathsf{A}}\).

By a similar argument we also have that
\[
\mathrm{tg}^{k\mathsf{A}}(g\#^{k\mathsf{A}}f)
=
\mathrm{tg}^{k\mathsf{A}}(g).
\]

\textsf{(C4)}
For every \(f\in A_{\omega}\), we define 
\begin{enumerate}
\item
\(n_{k,f}=\min\{i\in\mathbb{N} \mid k+1\leq i, f\in A_{i}\}\) and
\item
\(n_{k,f,\mathrm{sc}(f)}=\min\{i\in\mathbb{N} \mid k+1\leq i, f,\mathrm{sc}^{k\mathsf{A}_{n_{k,f}}}(f)\in A_{i}\}\).
\end{enumerate}
Let us note that \(n_{k,f}\leq n_{k,f,\mathrm{sc}(f)}\).

The following chain of equalities holds
\begin{align*}
f\#^{k\mathsf{A}}\mathrm{sc}^{k\mathsf{A}}(f)
&=
f\#^{k\mathsf{A}}\mathrm{sc}^{k\mathsf{A}_{n_{k,f}}}(f)
\tag{1}
\\
&=
f\#^{k\mathsf{A}_{n_{k,f,\mathrm{sc}(f)}}}\mathrm{sc}^{k\mathsf{A}_{n_{k,f}}}(f)
\tag{2}
\\
&=
f\#^{k\mathsf{A}_{n_{k,f,\mathrm{sc}(f)}}}\mathrm{sc}^{k\mathsf{A}_{n_{k,f,\mathrm{sc}(f)}}}(f)
\tag{3}
\\
&=f.
\tag{4}
\end{align*}

The first and second equalities unravel the definitions of \(\mathrm{sc}^{k\mathsf{A}}\) and \(\#^{k\mathsf{A}}\), respectively; the third equality follows from the fact that \(n_{k,f}\leq n_{k,f,\mathrm{sc}(f)}\); finally, the last equality follows from item~(C4) in Definition~\ref{DCat} since \((A_{n_{k,f,\mathrm{sc}(f)}},\xi^{\mathsf{A}_{n_{k,f,\mathrm{sc}(f)}}}_{k})\) is a single-sorted category.

By a similar argument we also have that
\[
\mathrm{tg}^{k\mathsf{A}}(f)\#^{k\mathsf{A}}f
=
f.
\]

\textsf{(C5)}
For every \(f,g,h\in A_{\omega}\), if \(\mathrm{sc}^{k\mathsf{A}}(h)=\mathrm{tg}^{k\mathsf{A}}(g)\) and \(\mathrm{sc}^{k\mathsf{A}}(g)=\mathrm{tg}^{k\mathsf{A}}(f)\), we define
\begin{enumerate}
\item
\(n_{k,f,g}=\min\{i\in\mathbb{N} \mid k+1\leq i, f,g\in A_{i}\}\),
\item
\(n_{k,h,g\#f}=\min\{i\in\mathbb{N} \mid k+1\leq i, h,g\#^{k\mathsf{A}_{n_{k,f,g}}}f\in A_{i}\}\),
\item
\(n_{k,g,h}=\min\{i\in\mathbb{N} \mid k+1\leq i, g,h\in A_{i}\}\),
\item
\(n_{k,h\#g,f}=\min\{i\in\mathbb{N} \mid k+1\leq i h\#^{k\mathsf{A}_{n_{k,g,h}}}g,f\in A_{i}\}\), and
\item
\(n_{k,f,g,h}=\min\{i\in\mathbb{N} \mid k+1\leq i, f,g,h\in A_{i}\}\).
\end{enumerate}
Let us note that \(n_{k,h,g\#f}\leq n_{k,f,g,h}\), \(n_{k,f,g}\leq n_{k,f,g,h}\), \(n_{k,h\#g,f}\leq n_{k,f,g,h}\) and \(n_{k,g,h}\leq n_{k,f,g,h}\).

The following chain of equalities holds
\allowdisplaybreaks
\begin{align*}
h\#^{k\mathsf{A}}(g\#^{k\mathsf{A}}f)
&=
h\#^{k\mathsf{A}}(g\#^{k\mathsf{A}_{n_{k,f,g}}}f)
\tag{1}
\\
&=
h\#^{k\mathsf{A}_{n_{k,h,g\#f}}}(g\#^{k\mathsf{A}_{n_{k,f,g}}}f)
\tag{2}
\\
&=
h\#^{k\mathsf{A}_{n_{k,f,g,h}}}(g\#^{k\mathsf{A}_{n_{k,f,g,h}}}f)
\tag{3}
\\
&=
(h\#^{k\mathsf{A}_{n_{k,f,g,h}}}g)\#^{k\mathsf{A}_{n_{k,f,g,h}}}f
\tag{4}
\\
&=
(h\#^{k\mathsf{A}_{n_{k,g,h}}}g)\#^{k\mathsf{A}_{n_{k,h\#g,f}}}f
\tag{5}
\\
&=
(h\#^{k\mathsf{A}_{n_{k,g,h}}}g)\#^{k\mathsf{A}}f
\tag{6}
\\
&=
(h\#^{k\mathsf{A}}g)\#^{k\mathsf{A}}f.
\tag{7}
\end{align*}

The first and second equality unravel the definition of \(\#^{k\mathsf{A}}\); the third equality follows from the fact that \(n_{k,h,g\#f}\leq n_{k,f,g,h}\) and \(n_{k,f,g}\leq n_{k,f,g,h}\); the fourth equality follows from item~(C5) in Definition~\ref{DCat} since \((A_{n_{k,f,g,h}},\xi^{\mathsf{A}_{n_{k,f,g,h}}}_{k})\) is a single-sorted category; the fifth equality follows from the fact that \(n_{k,h\#g,f}\leq n_{k,f,g,h}\) and \(n_{k,g,h}\leq n_{k,f,g,h}\); finally, the last two equalities recover the definition of \(\#^{k\mathsf{A}}\).

It follows that, for every \(k\in\omega\), \((A_{\omega},\xi^{\mathsf{A}}_{k})\) is a single-sorted category.

According to Definition~\ref{DnCatSS}, now we have to prove that, for every \(j,k\in\omega\) with \(j<k\), \((A_{\omega},(\xi^{\mathsf{A}}_{i})_{i\in\{j,k\}})\) is a single-sorted 2-category. To this end, we consider the different items stated in Definition~\ref{D2Cat}.

\textsf{(2C1)}
For every \(f\in A_{\omega}\), we define
\begin{enumerate}
\item
\(n_{k,f}=\min\{i\in\mathbb{N} \mid k+1\leq i, f\in A_{i}\}\),
\item
\(n_{j,\mathrm{sc}(f)}=\min\{i\in\mathbb{N} \mid j+1\leq i, \mathrm{sc}^{k\mathsf{A}_{n_{k,f}}}(f)\in A_{i}\}\), and
\item
\(n_{j,f}=\min\{i\in\mathbb{N} \mid j+1\leq i, f\in A_{i}\}\).
\end{enumerate}
Let us note that \(n_{j,\mathrm{sc}(f)}\leq n_{k,f}\) and \(n_{j,f}\leq n_{k,f}\).

The following chain of equalities holds
\begin{align*}
\mathrm{sc}^{j\mathsf{A}}(\mathrm{sc}^{k\mathsf{A}}(f))
&=
\mathrm{sc}^{j\mathsf{A}}(\mathrm{sc}^{k\mathsf{A}_{n_{k,f}}}(f))
\tag{1}
\\
&=
\mathrm{sc}^{j\mathsf{A}_{n_{j,\mathrm{sc}(f)}}}(\mathrm{sc}^{k\mathsf{A}_{n_{k,f}}}(f))
\tag{2}
\\
&=
\mathrm{sc}^{j\mathsf{A}_{n_{k,f}}}(\mathrm{sc}^{k\mathsf{A}_{n_{k,f}}}(f))
\tag{3}
\\
&=
\mathrm{sc}^{j\mathsf{A}_{n_{k,f}}}(f)
\tag{4}
\\
&=
\mathrm{sc}^{j\mathsf{A}_{n_{j,f}}}(f)
\tag{5}
\\
&=
\mathrm{sc}^{j\mathsf{A}}(f).
\tag{6}
\end{align*}

The first and second equality unravel the definition of \(\mathrm{sc}^{k\mathsf{A}}\) and \(\mathrm{sc}^{j\mathsf{A}}\) respectively; the third equality follows from the fact that \(n_{j,\mathrm{sc}(f)}\leq n_{k,f}\), the fourth equality follows item~(2C1) in Definition~\ref{D2Cat} since \((A_{n_{k,f}},(\xi^{\mathsf{A}_{n_{k,f}}}_{i})_{i\in\{j,k\}})\) is a 2-category; the fifth equality follows from the fact that \(n_{j,f}\leq n_{k,f}\); finally, the last equality recovers the definition of \(\mathrm{sc}^{j\mathsf{A}}\)

By a similar argument we also have that
\[
\mathrm{sc}^{j\mathsf{A}}\circ\mathrm{tg}^{k\mathsf{A}}=\mathrm{sc}^{j\mathsf{A}}.
\]

On the other hand, for every \(f\in A_{\omega}\), we define
\begin{enumerate}
\item
\(n_{j,f}=\min\{i\in\mathbb{N} \mid j+1\leq i, f\in A_{i}\}\),
\item
\(n_{k,\mathrm{sc}(f)}=\min\{i\in\mathbb{N} \mid k+1\leq i, \mathrm{sc}^{j\mathsf{A}_{n_{j,f}}}(f)\in A_{i}\}\) and
\item
\(n_{k,f}=\min\{i\in\mathbb{N} \mid k+1\leq i, f\in A_{i}\}\).
\end{enumerate}
Let us note that \(n_{j,f}\leq n_{k,f}\) and \(n_{k,\mathrm{sc}(f)}\leq n_{k,f}\).

The following chain of equalities holds
\allowdisplaybreaks
\begin{align*}
\mathrm{sc}^{k\mathsf{A}}(\mathrm{sc}^{j\mathsf{A}}(f))
&=
\mathrm{sc}^{k\mathsf{A}}(\mathrm{sc}^{j\mathsf{A}_{n_{j,f}}}(f))
\tag{1}
\\
&=
\mathrm{sc}^{k\mathsf{A}_{n_{k,\mathrm{sc}(f)}}}(\mathrm{sc}^{j\mathsf{A}_{n_{j,f}}}(f))
\tag{2}
\\
&=
\mathrm{sc}^{k\mathsf{A}_{n_{k,f}}}(\mathrm{sc}^{j\mathsf{A}_{n_{k,f}}}(f))
\tag{3}
\\
&=
\mathrm{sc}^{j\mathsf{A}_{n_{k,f}}}(f)
\tag{4}
\\
&=
\mathrm{sc}^{j\mathsf{A}_{n_{j,f}}}(f)
\tag{5}
\\
&=
\mathrm{sc}^{j\mathsf{A}}(f).
\tag{6}
\end{align*}

The first and second equality unravel the definition of \(\mathrm{sc}^{j\mathsf{A}}\) and \(\mathrm{sc}^{k\mathsf{A}}\), respectively; the third equality follows from the fact that \(n_{j,f}\leq n_{k,f}\) and \(n_{k,\mathrm{sc}(f)}\leq n_{k,f}\); the fourth equality follows from item~(2C1) in Definition~\ref{D2Cat} since \((A_{n_{k,f}},(\xi^{\mathsf{A}_{n_{k,f}}}_{i})_{i\in\{j,k\}})\) is a \(2\)-category; the fifth equality follows from the fact that \(n_{j,f}\leq n_{k,f}\); finally the last equality recovers the definition of \(\mathrm{sc}^{j\mathsf{A}}\).

By a similar argument we also have that
\[
\mathrm{tg}^{k\mathsf{A}}\circ\mathrm{tg}^{j\mathsf{A}}
=
\mathrm{tg}^{j\mathsf{A}}
=
\mathrm{tg}^{j\mathsf{A}}\circ\mathrm{tg}^{k\mathsf{A}}
=
\mathrm{tg}^{j\mathsf{A}}\circ\mathrm{sc}^{k\mathsf{A}}.
\]

\textsf{(2C2)}
For every \(k,j\in A_{\omega}\), if \(\mathrm{sc}^{j\mathsf{A}}(g)=\mathrm{tg}^{k\mathsf{A}}(f)\), we define
\begin{enumerate}
\item
\(n_{j,f,g}=\min\{i\in\mathbb{N} \mid j+1\leq i, f,g\in A_{i}\}\),
\item
\(n_{k,g\#f}=\min\{i\in\mathbb{N} \mid k+1\leq i, g\#^{j\mathsf{A}_{n_{j,f,g}}}f\in A_{i}\}\),
\item
\(n_{k,f}=\min\{i\in\mathbb{N} \mid k+1\leq i, f\in A_{i}\}\),
\item
\(n_{k,g}=\min\{i\in\mathbb{N} \mid k+1\leq i, g\in A_{i}\}\),
\item
\(n_{k,f,g}=\min\{i\in\mathbb{N} \mid k+1\leq i, f,g\in A_{i}\}\), and 
\item
\(n_{j,\mathrm{sc}(g),\mathrm{sc}(f)}=\min\{i\in\mathbb{N} \mid j+1\leq i, \mathrm{sc}^{k\mathsf{A}_{n_{k,g}}}(g),\mathrm{sc}^{k\mathsf{A}_{n_{k,f}}}(f)\in A_{i}\}\).
\end{enumerate}
Let us note that \(n_{k,g\#f}\leq n_{k,f,g}\), \(n_{j,f,g}\leq n_{k,f,g}\), \(n_{k,f}\leq n_{k,f,g}\), \(n_{k,g}\leq n_{k,f,g}\) and \(n_{j,\mathrm{sc}(g),\mathrm{sc}(f)}\leq n_{k,f,g}\).

The following chain of equalities holds
\begin{align*}
\mathrm{sc}^{k\mathsf{A}}(g\#^{j\mathsf{A}}f)
&=
\mathrm{sc}^{k\mathsf{A}}(g\#^{j\mathsf{A}_{n_{j,f,g}}}f)
\tag{1}
\\
&=
\mathrm{sc}^{k\mathsf{A}_{n_{k,g\#f}}}(g\#^{j\mathsf{A}_{n_{j,f,g}}}f)
\tag{2}
\\
&=
\mathrm{sc}^{k\mathsf{A}_{n_{k,f,g}}}(g\#^{j\mathsf{A}_{n_{k,f,g}}}f)
\tag{3}
\\
&=
\mathrm{sc}^{k\mathsf{A}_{n_{k,f,g}}}(g)\#^{j\mathsf{A}_{n_{k,f,g}}}\mathrm{sc}^{k\mathsf{A}_{n_{k,f,g}}}(f)
\tag{4}
\\
&=
\mathrm{sc}^{k\mathsf{A}_{n_{k,g}}}(g)\#^{j\mathsf{A}_{n_{j,\mathrm{sc}(g),\mathrm{sc}(f)}}}\mathrm{sc}^{k\mathsf{A}_{n_{k,f}}}(f)
\tag{5}
\\
&=
\mathrm{sc}^{k\mathsf{A}_{n_{k,g}}}(g)\#^{j\mathsf{A}}\mathrm{sc}^{k\mathsf{A}_{n_{k,f}}}(f)
\tag{6}
\\
&=
\mathrm{sc}^{k\mathsf{A}}(g)\#^{j\mathsf{A}}\mathrm{sc}^{k\mathsf{A}}(f).
\tag{7}
\end{align*}

The first and second equality unravel the definition of \(\#^{j\mathsf{A}}\) and \(\mathrm{sc}^{k\mathsf{A}}\) respectively; the third equality follows from the fact that \(n_{k,g\#f}\leq n_{k,f,g}\) and \(n_{j,f,g}\leq n_{k,f,g}\); the fourth equality follows from item~(2C2) in Definition~\ref{D2Cat} since \((A_{n'},(\xi^{\mathsf{A}_{n'}}_{i})_{i\in\{j,k\}})\) is a \(2\)-category; the fifth equality follows from the fact that \(n_{k,f}\leq n_{k,f,g}\), \(n_{k,g}\leq n_{k,f,g}\) and \(n_{j,\mathrm{sc}(g),\mathrm{sc}(f)}\leq n_{k,f,g}\); the last two equalities recover the definition of \(\#^{j\mathsf{A}}\) and \(\mathrm{sc}^{k\mathsf{A}}\), respectively.

By a similar argument we also have that
\[
\mathrm{tg}^{k\mathsf{A}}(g\#^{j\mathsf{A}}f)=\mathrm{tg}^{k\mathsf{A}}(g)\#^{j\mathsf{A}}\mathrm{tg}^{k\mathsf{A}}(f).
\]

\textsf{(2C3)}
For every \(f,f',g,g'\in A_{\omega}\) if
\begin{align*}
\mathrm{sc}^{k\mathsf{A}}(g')&=\mathrm{tg}^{k\mathsf{A}}(g);
&\mathrm{sc}^{k\mathsf{A}}(f')&=\mathrm{tg}^{k\mathsf{A}}(f);
\\
\mathrm{sc}^{j\mathsf{A}}(g')&=\mathrm{tg}^{j\mathsf{A}}(f');
&\mathrm{sc}^{j\mathsf{A}}(g)&=\mathrm{tg}^{j\mathsf{A}}(f),
\end{align*}
we define
\begin{enumerate}
\item
\(n_{j,f,g}=\min\{i\in\mathbb{N} \mid j+1\leq i, f,g\in A_{i}\}\),
\item
\(n_{j,f',g'}=\min\{i\in\mathbb{N} \mid j+1\leq i, f',g'\in A_{i}\}\),
\item
\(n_{j,f,f'}=\min\{i\in\mathbb{N} \mid j+1\leq i, f,f'\in A_{i}\}\),
\item
\(n_{j,g,g'}=\min\{i\in\mathbb{N} \mid j+1\leq i, g,g'\in A_{i}\}\),
\item
\(n_{k,g\#f,g'\#f'}=\min\{i\in\mathbb{N} \mid k+1\leq i, g\#^{j\mathsf{A}_{n_{j,f,g}}}f,g'\#^{j\mathsf{A}_{n_{j,f',g'}}}f'\in A_{i}\}\),
\item
\(n_{k,g'\#g,f'\#f}=\min\{i\in\mathbb{N} \mid k+1\leq i, g'\#^{j\mathsf{A}_{n_{j,g,g'}}}g,f'\#^{j\mathsf{A}_{n_{j,f,f'}}}f\in A_{i}\}\), and
\item
\(n_{k,f,f',g,g'}=\min\{i\in\mathbb{N} \mid k+1\leq i, f,f',g,g'\in A_{i}\}\).
\end{enumerate}
Let us note that \(n_{j,f,g}\leq n_{k,f,f',g,g'}\), \(n_{j,f',g'}\leq n_{k,f,f',g,g'}\), \(n_{k,g\#f,g'\#f'}\leq n_{k,f,f',g,g'}\), \(n_{j,f,f'}\leq n_{k,f,f',g,g'}\), \(n_{j,g,g'}\leq n_{k,f,f',g,g'}\) and \(n_{k,g'\#g,f'\#f}\leq n_{k,f,f',g,g'}\).

The following chain of equalities holds
\begin{align*}
&(g'\#^{j\mathsf{A}}f')\#^{k\mathsf{A}}(g\#^{j\mathsf{A}}f)
\\
&=
(g'\#^{j\mathsf{A}_{n_{j,f',g'}}}f')\#^{k\mathsf{A}}(g\#^{j\mathsf{A}_{n_{j,f,g}}}f)
\tag{1}
\\
&=
(g'\#^{j\mathsf{A}_{n_{j,f',g'}}}f')\#^{k\mathsf{A}_{n_{k,g\#f,g'\#f'}}}(g\#^{j\mathsf{A}_{n_{j,f,g}}}f)
\tag{2}
\\
&=
(g'\#^{j\mathsf{A}_{n_{k,f,f',g,g'}}}f')\#^{k\mathsf{A}_{n_{k,f,f',g,g'}}}(g\#^{j\mathsf{A}_{n_{k,f,f',g,g'}}}f)
\tag{3}
\\
&=
(g'\#^{j\mathsf{A}_{n_{k,f,f',g,g'}}}g)\#^{k\mathsf{A}_{n_{k,f,f',g,g'}}}(f'\#^{j\mathsf{A}_{n_{k,f,f',g,g'}}}f)
\tag{4}
\\
&=
(g'\#^{j\mathsf{A}_{n_{j,g,g'}}}g)\#^{k\mathsf{A}_{n_{k,g'\#g,f'\#f}}}(f'\#^{j\mathsf{A}_{n_{j,f,f'}}}f)
\tag{5}
\\
&=
(g'\#^{j\mathsf{A}_{n_{j,g,g'}}}g)\#^{k\mathsf{A}}(f'\#^{j\mathsf{A}_{n_{j,f,f'}}}f)
\tag{6}
\\
&=
(g'\#^{j\mathsf{A}}g)\#^{k\mathsf{A}}(f'\#^{j\mathsf{A}}f).
\tag{7}
\end{align*}

The first and second equalities unravel the definition of \(\#^{j\mathsf{A}}\) and \(\#^{k\mathsf{A}}\), respectively; the third equality holds because \(n_{j,f,g}\leq n_{k,f,f',g,g'}\), \(n_{j,f',g'}\leq n_{k,f,f',g,g'}\) and \(n_{k,g\#f,g'\#f'}\leq n_{k,f,f',g,g'}\); the fourth equality follows from item~(2C3) in Definition~\ref{D2Cat} since \((A_{n_{k,f,f',g,g'}},(\xi^{\mathsf{A}_{n_{k,f,f',g,g'}}}_{i})_{i\in \{j,k\}})\) is a \(2\)-category; the fifth equality follows from the fact that \(n_{j,f,f'}\leq n_{k,f,f',g,g'}\), \(n_{j,g,g'}\leq n_{k,f,f',g,g'}\) and \(n_{k,g'\#g,f'\#f}\leq n_{k,f,f',g,g'}\); the last two equalities recover the definition of \(\#^{k\mathsf{A}}\) and \(\#^{j\mathsf{A}}\), respectively.

It follows that, for every \(j,k\in\omega\) with \(j<k\), \((A_{\omega}, (\xi^{\mathsf{A}}_{i})_{i\in\{j,k\}})\) is a single-sorted \(2\)-category.

Finally, following Definition~\ref{DnCatSS}, we have to prove that, for every \(f\in A_{\omega}\), there exists \(k\in\omega\) such that \(f=\mathrm{sc}^{k\mathsf{A}}(f)\) or, equivalently, that \(A_{\omega}=\bigcup_{k\in\omega}(A_{\omega})_{k}\).

Let \(f\in A_{\omega}\) and \(k\in\omega\). We prove that
\begin{align*}
f&=\mathrm{sc}^{k\mathsf{A}}(f)
&\mbox{if, and only if,}&&
f\in A_{k}.
\end{align*}

Let us suppose that \(f=\mathrm{sc}^{k\mathsf{A}}(f)\). Since \(f\in\bigcup_{k\in\omega}A_{k}\), let \(n_{f}\) be the minimum natural number \(i\in\mathbb{N}\) for which \(f\in A_{i}\), that is, \(n_{f}=\min\{i\in\mathbb{N} \mid f\in A_{i}\}\). We distinguish two cases (1) \(n_{f}\leq k\); and (2) \(n_{f}>k\). If (1), then \(f\in A_{n_{f}}\subseteq A_{k}\) where the last inclusion follows from item~(CC1). If (2), then \(n_{f}=\min\{i\in\mathbb{N} \mid k+1\leq i, f\in A_{i}\}\). Therefore, by the definition of \(\mathrm{sc}^{k\mathsf{A}}\),
\[
f
=
\mathrm{sc}^{k\mathsf{A}}(f)
=
\mathrm{sc}^{k\mathsf{A}_{n_{f}}}(f).
\]
That implies that \(f\) is a \(k\)-cell in \((A_{n_{f}},\xi^{\mathsf{A}_{n_{f}}}_{k})\) and, by item~(CC1), \(f\in (A_{n_{f}})_{k}=A_{k}\) contradicting the minimality of \(n_{f}\).

Reciprocally, let us suppose that \(f\in A_{k}\). Then the following chain of  hold
\begin{align*}
\mathrm{sc}^{k\mathsf{A}}(f)
&=
\mathrm{sc}^{k\mathsf{A}_{k+1}}(f)
\tag{1}
\\
&=
f.
\tag{2}
\end{align*}

The first equality follows from the definition of \(\mathrm{sc}^{k\mathsf{A}}\), since 
\[
\min\{i\in\mathbb{N} \mid k+1\leq i, f\in A_{i}\}=k+1;
\]
and the second equality follows from the fact that, by item~(CC1), \(f\in A_{k}=(A_{k+1})_{k}\).

All in all, we can affirm that \(\mathsf{A}_{\omega}\) is a single-sorted \(\omega\)-category.

We now want to prove that, for \(n\in\omega\), \(\mathsf{A}_{\omega}^{<n}=\mathsf{A}_{n}\).

The fact that \((A_{\omega})_{n}=A_{n}\) proves that the underlying set of \(\mathsf{A}_{\omega}^{<n}\) is equal to the underlying set of \(\mathsf{A}_{n}\). Let \(k\in n\) and \(f\in A_{n}\). Then the following chain of equalities holds
\begin{align*}
\mathrm{sc}^{k\mathsf{A}}\bigr|^{A_{n}}_{A_{n}}(f)
&=
\mathrm{sc}^{k\mathsf{A}}(f)
\tag{1}
\\
&=
\mathrm{sc}^{k\mathsf{A}_{n}}(f)
\tag{2}
\end{align*}

The first equality follows from the fact that \(f\in A_{n}\); the second equality follows from the definition of \(\mathrm{sc}^{k\mathsf{A}}\), because 
\[
\min\{i\in\mathbb{N}\mid k+1\leq i, f\in A_{i}\}\leq n.
\]

By a similar argument we also have that
\begin{align*}
\mathrm{tg}^{k\mathsf{A}}\bigr|^{A_{n}}_{A_{n}}=\mathrm{tg}^{k\mathsf{A}_{n}}&
&\mbox{and}&&
\#^{k\mathsf{A}}\bigr|^{A_{n}}_{A_{n}}=\#^{k\mathsf{A}_{n}}
\end{align*}

It follows that, for every \(k\in n\), \((\xi^{\mathsf{A}}_{k})^{<n}=\xi^{\mathsf{A}_{n}}_{k}\).

All in all, we can affirm that \(\mathsf{A}_{\omega}^{<n}=\mathsf{A}_{n}\).

Finally, we prove that \(\mathsf{A}_{\omega}\) is the unique construction with such property.

Let \(\mathsf{C}=(C,(\xi_{k})_{k\in\omega})\) be a single-sorted \(\omega\)-category with \((\xi_{k})_{k\in\omega}=(\#^{k},\mathrm{sc}^{k},\mathrm{tg}^{k})_{k\in\omega}\) such that, for every \(n\in\omega\),
\[
\mathsf{C}^{<n}=\mathsf{A}_{n}.
\]

We want to prove that \(\mathsf{C}=\mathsf{A}_{\omega}\), that is, \(C=A_{\omega}\) and, for every \(k\in\omega\), \(\xi_{k}=\xi^{\mathsf{A}}_{k}\).

By assumption, for every \(n\in\omega\), \(\mathsf{C}^{<n}=\mathsf{A}_{n}\), therefore \(C_{n}=A_{n}\). Then the following chain of equalities holds
\begin{align*}
C
&=
\textstyle\bigcup_{k\in\omega}C_{k}
\tag{1}
\\
&=
\textstyle\bigcup_{k\in\omega}A_{k}
\tag{2}
\\
&=
A_{\omega}.
\tag{3}
\end{align*}

The first equality follows from item~(\(\omega\)C2) of the definition of single-sorted \(\omega\)-category in Definition~\ref{DnCatSS}; the second equality follows from the fact that \(C_{k}=A_{k}\); finally, the last equality recovers the definition of \(A_{\omega}\).

Let \(k\in\omega\) and \(f\in A_{\omega}\). Let 
\[
n_{k,f}=\min\{i\in\mathbb{N} \mid k+1\leq i, f\in A_{i}\}.
\]
By assumption, \(\mathsf{C}^{<n_{k,f}}=\mathsf{A}_{n_{k,f}}\), therefore \((\xi_{k})_{k\in n_{k,f}}=(\xi_{k}^{A_{n_{k,f}}})_{k\in n_{k,f}}\). In particular, since \(k<n_{k,f}\), \(\xi_{k}=\xi^{A_{n_{k,f}}}_{k}\). Then the following chain of equalities holds
\begin{align*}
\mathrm{sc}^{k\mathsf{A}}(f)
&=
\mathrm{sc}^{k\mathsf{A}_{n_{k,f}}}(f)
\tag{1}
\\
&=
\mathrm{sc}^{k}(f).
\tag{2}
\end{align*}

The first equality unravels the definition of \(\mathrm{sc}^{k\mathsf{A}}\); the second equality follows from the fact that \(\xi_{k}=\xi^{\mathsf{A}_{n_{k,f}}}_{k}\).

By a similar argument we also have that, for every \(f,g\in A_{\omega}\),
\begin{align*}
\mathrm{tg}^{n\mathsf{A}}(f)=\mathrm{tg}^{n}(f)&
&\mbox{and}&&
f\#^{n\mathsf{A}}g=f\#^{n}g.
\end{align*}

This proves that \(\mathsf{C}=\mathsf{A}_{\omega}\).

Claim~\ref{CUnivPropSSA} follows.

\begin{claim}\label{CUnivPropSSB}
For every morphism \(\varphi\) from \(A\) to \(B\) in \(\mathsf{L}\) we can uniquely assign a single-sorted \(\omega\)-functor \(\varphi_{\omega}\) from \(\mathsf{A}_{\omega}\) to \(\mathsf{B}_{\omega}\) such that, for every \(n\in\omega\),
\[
\varphi_{\omega}^{<n}=\varphi_{n}.
\]
\end{claim}

Let us point out that, for every \(m,n\in\omega\) with \(m\leq n\), since \(F_{m}=U^{(m,n)}\circ F_{n}\), it follows that \(\varphi_{m}=\varphi_{n}^{<m}\). In particular,
\begin{enumerate}
\item[(CC3)]
For every \(f\in A_{m}\), \(\varphi_{m}(f)=\varphi_{n}(f)\bigr|^{B_{m}}_{A_{m}}\).
\end{enumerate}

Let \(\varphi\) be a morphism from \(A\) to \(B\) in \(\mathsf{L}\). We consider the mapping \(\varphi_{\omega}\) from \(A_{\omega}\) to \(B_{\omega}\) defined, for every \(f\in A_{\omega}\), as \(\varphi_{n_{f}}(f)\) where
\[
n_{f}=\min\{i\in\mathbb{N} \mid f\in A_{i}\}.
\]
Let us note that, from item~(CC3) above, it follows that, for every \(n\in\mathbb{N}\) with \(n_{f}\geq n\), \(\varphi_{\omega}(f)=\varphi_{n}(f)\). We want to prove that \(\varphi_{\omega}\) is a \(\omega\)-functor. Thus, we must check that it satisfies the conditions stated in Definition~\ref{DnCatSS}.

We must prove that, for \(k\in\omega\), \(\varphi_{\omega}\) is a single-sorted functor from the single-sorted category \((A_{\omega},\xi^{\mathsf{A}}_{k})\) to the single-sorted category \((B_{\omega},\xi^{\mathsf{B}}_{k})\). To this end, we consider the different items stated in Definition~\ref{DCat}.

\textsf{(Ci)}
For every \(f\in A_{\omega}\), we define
\begin{enumerate}
\item
\(n_{k,f}=\min\{i\in\mathbb{N} \mid k+1\leq i, f\in A_{i}\}\),
\item
\(n_{\mathrm{sc}(f)}=\min\{i\in\mathbb{N} \mid \mathrm{sc}^{k\mathsf{A}_{n_{k,f}}}(f)\in A_{i}\}\),
\item
\(n_{k,\varphi(f)}=\min\{i\in\mathbb{N} \mid k+1\leq i, \varphi_{\omega}(f)\in B_{i}\}\), and
\item
\(n_{f}=\min\{i\in\mathbb{N} \mid f\in A_{i}\}\).
\end{enumerate}
Let us note that \(n_{\mathrm{sc}(f)}\leq n_{k,f}\), \(n_{k,\varphi(f)}\leq n_{k,f}\) and \(n_{f}\leq n_{k,f}\).

The following chain of equalities holds
\begin{align*}
\varphi_{\omega}(\mathrm{sc}^{k\mathsf{A}}(f))
&=
\varphi_{\omega}(\mathrm{sc}^{k\mathsf{A}_{n_{k,f}}}(f))
\tag{1}
\\
&=
\varphi_{n_{\mathrm{sc}(f)}}(\mathrm{sc}^{k\mathsf{A}_{n_{k,f}}}(f))
\tag{2}
\\
&=
\varphi_{n}(\mathrm{sc}^{k\mathsf{A}_{n_{k,f}}}(f))
\tag{3}
\\
&=
\mathrm{sc}^{k\mathsf{B}_{n_{k,f}}}(\varphi_{n_{k,f}}(f))
\tag{4}
\\
&=
\mathrm{sc}^{k\mathsf{B}_{n_{k,\varphi(f)}}}(\varphi_{n_{f}}(f))
\tag{5}
\\
&=
\mathrm{sc}^{k\mathsf{B}_{n_{k,\varphi(f)}}}(\varphi_{\omega}(f))
\tag{6}
\\
&=
\mathrm{sc}^{k\mathsf{B}}(\varphi_{\omega}(f)).
\tag{7}
\end{align*}

The first and second equality unravels the definition of \(\mathrm{sc}^{k\mathsf{A}}\) and \(\varphi_{\omega}\), respectively; the third equality follows from the fact that \(n_{\mathrm{sc}(f)}\leq n_{k,f}\); the fourth equality follows from item~(Ci) in Definition~\ref{DCat} since \(\varphi_{n_{k,f}}\) is a functor from \((A_{n_{k,f}},\xi^{\mathsf{A}_{n_{k,f}}}_{k})\) to \((B_{n_{k,f}},\xi^{\mathsf{B}_{n_{k,f}}}_{k})\); the fifth equality follows from the fact that \(n_{k,\varphi(f)}\leq n_{k,f}\) and \(n_{f}\leq n_{k,f}\); the last two equalities recover the definition of \(\varphi_{\omega}\) and \(\mathrm{sc}^{k\mathsf{B}}\), respectively.

By a similar argument we also have that 
\[
\varphi_{\omega}\circ\mathrm{tg}^{k\mathsf{A}}
=
\mathrm{tg}^{k\mathsf{B}}\circ\varphi_{\omega}.
\]

\textsf{(Cii)}
For every \(f,g\in A_{\omega}\), if \(\mathrm{sc}^{k\mathsf{A}}(g)=\mathrm{tg}^{k\mathsf{A}}(f)\), then we define
\begin{enumerate}
\item
\(n_{k,f,g}=\min\{i\in\mathbb{N} \mid k+1\leq i, f,g\in A_{i}\}\),
\item
\(n_{g\#f}=\min\{i\in\mathbb{N}\mid g\#^{k\mathsf{A}_{n_{k,f,g}}}f\in A_{i}\}\),
\item
\(n_{f}=\min\{i\in\mathbb{N} \mid f\in A_{i}\}\),
\item
\(n_{g}=\min\{i\in\mathbb{N}\mid g\in A_{i}\}\), and
\item
\(n_{k,\varphi(f),\varphi(g)}=\min\{i\in\mathbb{N} \mid k+1\leq i, \varphi_{\omega}(f),\varphi_{\omega}(g)\in B_{i}\}\).
\end{enumerate}
Let us note that \(n_{g\#f}\leq n_{k,f,g}\), \(n_{f}\leq n_{k,f,g}\), \(n_{g}\leq n_{k,f,g}\) and \(n_{k,\varphi(f),\varphi(g)}\leq n_{k,f,g}\).

The following chain of equalities holds
\begin{align*}
\varphi_{\omega}(g\#^{k\mathsf{A}}f)
&=
\varphi_{\omega}(g\#^{k\mathsf{A}_{n_{k,f,g}}}f)
\tag{1}
\\
&=
\varphi_{n_{g\#f}}(g\#^{k\mathsf{A}_{n_{k,f,g}}}f)
\tag{2}
\\
&=
\varphi_{n_{k,f,g}}(g\#^{k\mathsf{A}_{n_{k,f,g}}}f)
\tag{3}
\\
&=
\varphi_{n_{k,f,g}}(g)\#^{k\mathsf{B}_{n_{k,f,g}}}\varphi_{n_{k,f,g}}(f)
\tag{4}
\\
&=
\varphi_{n_{g}}(g)\#^{k\mathsf{B}_{n_{k,\varphi(f),\varphi(g)}}}\varphi_{n_{f}}(f)
\tag{5}
\\
&=
\varphi_{\omega}(g)\#^{k\mathsf{B}_{n_{k,\varphi(f),\varphi(g)}}}\varphi_{\omega}(f)
\tag{6}
\\
&=
\varphi_{\omega}(g)\#^{k\mathsf{B}}\varphi_{\omega}(f).
\tag{7}
\end{align*}

The first and second equalities unravel the definition of \(\#^{k\mathsf{A}}\) and \(\varphi_{\omega}\), respectively; the third equality follows from the fact that \(n_{g\#f}\leq n_{k,f,g}\); the fourth equality follows from item~(Cii) in Definition~\ref{DCat} since \(\varphi_{n}\) is a functor from \((A_{n_{k,f,g}},\xi^{\mathsf{A}_{n_{k,f,g}}}_{k})\) to \((B_{n_{k,f,g}},\xi^{\mathsf{B}_{n_{k,f,g}}}_{k})\); the fifth equality follows from the fact that \(n_{f}\leq n_{k,f,g}\), \(n_{g}\leq n_{k,f,g}\) and \(n_{k,\varphi(f),\varphi(g)}\leq n_{k,f,g}\); the last two equalities recover the definition of \(\varphi_{\omega}\) and \(\#^{k\mathsf{B}}\), respectively.

All in all, we can affirm that \(\varphi_{\omega}\) is a \(\omega\)-functor. 

We now want to prove that, for \(n\in\omega\), \(\varphi_{\omega}^{<n}=\varphi_{n}\).

For every \(f\in A_{n}\), we define
\begin{enumerate}
\item
\(n_{f}=\min\{i\in\mathbb{N} \mid f\in A_{i}\}\).
\end{enumerate}

The following chain of equalities holds
\begin{align*}
\varphi_{\omega}^{<n}(f)
&=
\varphi_{\omega}\bigr|^{B_{n}}_{A_{n}}(f)
\tag{1}
\\
&=
\varphi_{\omega}(f)
\tag{2}
\\
&=
\varphi_{n_{f}}(f)
\tag{3}
\\
&=
\varphi_{n}(f).
\tag{4}
\end{align*}

The first equality unravels the definition of \(\varphi_{\omega}^{<n}\); the second equality follows from the fact that \(f\in A_{n}\); the third equality unravels the definition of \(\varphi_{\omega}\); finally, the last equality follows from the fact that \(n_{f}\leq n\).

All in all, we can affirm that \(\varphi_{\omega}^{<n}=\varphi_{n}\).

Finally, we prove that \(\varphi_{\omega}\) is the unique construction with such property.

Let \(\psi\) be a \(\omega\)-functor from 	\(\mathsf{A}_{\omega}\) to \(\mathsf{B}_{\omega}\), such that, for every \(n\in\omega\),
\[
\psi^{<n}=\varphi_{n}.
\]

We want to prove that \(\psi=\varphi_{\omega}\), that is, for every \(f\in A_{\omega}\), \(\psi(f)=\varphi_{\omega}(f)\).

For every \(f\in A_{\omega}\), defining
\begin{enumerate}
\item
\(n_{f}=\min\{i\in\mathbb{N} \mid f\in A_{i}\}\),
\end{enumerate}
the following chain of equalities holds
\begin{align*}
\varphi_{\omega}(f)
&=
\varphi_{n_{f}}(f)
\tag{1}
\\
&=
\psi^{<n_{f}}(f)
\tag{2}
\\
&=
\psi\bigr|^{B_{n_{f}}}_{A_{n_{f}}}(f)
\tag{3}
\\
&=
\psi(f).
\tag{4}
\end{align*}

The first equality unravels the definition of \(\varphi_{\omega}\); the second equality follows from the assumption that, for every \(n\in\omega\), \(\psi^{<n}=\varphi_{n}\); the third equality unravels the definition of \(\psi^{<n_{f}}\); the fourth equality follows from the fact that \(f\in A_{n_{f}}\).

Claim~\ref{CUnivPropSSB} follows. 

Therefore, we define the functor \(\langle F_{k}\rangle_{k\in\omega}\) from \(\mathsf{L}\) to \(\w\mathsf{Cat}\), whose object mapping sends an object \(A\) in \(\mathsf{L}\) to the \(\omega\)-category \(\mathsf{A}_{\omega}\), and whose morphism mapping sends a morphism \(\varphi\) from \(A\) to \(B\) in \(\mathsf{L}\) to the \(\omega\)-functor \(\varphi_{\omega}\) from \(\mathsf{A}_{\omega}\) to \(\mathsf{B}_{\omega}\).

It is straightforward to check that \(\langle F_{k}\rangle_{k\in\omega}\) preserves compositions and identities.

Claims~\ref{CUnivPropSSA} and \ref{CUnivPropSSB} prove that the functor is well-defined and that, for every \(n\in\omega\), the diagram
\[
\xymatrix{
\mathsf{L}
  \ar[d]_{\langle F_{k}\rangle_{k\in\omega}}
  \ar[dr]^{F_{n}}
\\
\w\mathsf{Cat}
  \ar[r]_{U^{(n,\omega)}}
  &
\mathsf{nCat}
}
\]
commutes. Moreover, from the unicity of both constructions follows the unicity of the functor \(\langle F_{k}\rangle_{k\in\omega}\).

This completes the proof.
\end{proof}

\section{Many-sorted presentation of higher-order categories}
We next define the notion of many-sorted $\omega$-category and in doing so we follow Leinster (see~\cite{L02}, p.~8). But before doing so, it is worth recalling that for the $\omega$-indexed family $(k)_{k\in \omega}$,  $\coprod_{k\in \omega}k$ is $\bigcup_{k\in \omega}(k\times \{k\})$. Thus, the elements of $\coprod_{k\in \omega}k$ are ordered pairs $(j,k)$ with $k\in \omega$ and $j\in k$. 

\begin{definition}\label{DnCatMS} 
A \emph{many-sorted $\omega$-category} is an ordered pair $\mathsf{C} = ((C_{k})_{k\in\omega},\zeta)$ consisting of an $\omega$-sorted set $(C_{k})_{k\in\omega}$ and a $\coprod_{k\in \omega}k$-indexed family
$$
\zeta = (\zeta_{j,k})_{(j,k)\in \coprod_{k\in \omega}k},
$$ 
where, for every $k\in \omega$ and every $j\in k$, $\zeta_{j,k} = (\circ^{(j)},\mathrm{sc}^{(j,k)},\mathrm{tg}^{(j,k)},\mathrm{i}^{(k,j)})$ with
\begin{enumerate}
\item $\mathrm{sc}^{(j,k)}$ and $\mathrm{tg}^{(j,k)}$ mappings from $C_{k}$ to $C_{j}$  
which send $f\in C_{k}$ to its \emph{$(j,k)$-source} and \emph{$(j,k)$-target}, respectively;

\item $\mathrm{i}^{(k,j)}$ a mapping from $C_{j}$ to $C_{k}$ 
which sends $f\in C_{j}$ to its \emph{$(k,j)$-identity}; and
%
\item $\circ^{(j)}$ a partial binary operation 
$$
\circ^{(j)}\colon C_{k}\times C_{k} \dmor C_{k},
$$ 
whose domain of definition is
$$
\mathrm{Dom}(\circ^{(j)}) = \left\lbrace(f,g)\in C_{k}\times C_{k}
\,
\middle|
\,
 \mathrm{sc}^{(j,k)}(g)=\mathrm{tg}^{(j,k)}(f)\right\rbrace,
$$ 
i.e., the fibered product of $C_{k}$ and $C_{k}$ over $\mathrm{sc}^{(j,k)}$ and $\mathrm{tg}^{(j,k)}$, and which sends a pair $(f,g)\in \mathrm{Dom}(\circ^{(j)})$ to its \emph{$j$-composition in $C_{k}$}.
\end{enumerate}
These data are required to satisfy the following conditions.
\begin{itemize}
\item[(MS1)]\label{MS1} For every $0\leq i<j<k$, if $f\in C_{k}$, then
\allowdisplaybreaks
\begin{align*}
\mathrm{sc}^{(i,j)}(\mathrm{sc}^{(j,k)}(f))&=\mathrm{sc}^{(i,k)}(f)=\mathrm{sc}^{(i,j)}(\mathrm{tg}^{(j,k)}(f)),\\
\mathrm{tg}^{(i,j)}(\mathrm{tg}^{(j,k)}(f))&=\mathrm{tg}^{(i,k)}(f)=\mathrm{tg}^{(i,j)}(\mathrm{sc}^{(j,k)}(f)).
\end{align*}
\item[(MS2)]\label{MS2} For every $0\leq j<k$, if $g,f\in C_{k}$ are such that $\mathrm{sc}^{(j,k)}(g)=\mathrm{tg}^{(j,k)}(f)$, then
\allowdisplaybreaks
\begin{align*}
\mathrm{sc}^{(j,k)}(g\circ^{(j)}f)&=\mathrm{sc}^{(j,k)}(f),
&
\mathrm{tg}^{(j,k)}(g\circ^{(j)}f)&=\mathrm{tg}^{(j,k)}(g).
\end{align*}
\item[(MS3)]\label{MS3} For every $0\leq i<j< k$, if $g,f\in C_{k}$ are such that $\mathrm{sc}^{(i,k)}(g)=\mathrm{tg}^{(i,k)}(f)$, then
\allowdisplaybreaks
\begin{align*}
\mathrm{sc}^{(j,k)}(g\circ^{(i)}f)&=\mathrm{sc}^{(j,k)}(g)\circ^{(i)}\mathrm{sc}^{(j,k)}(f),\\
\mathrm{tg}^{(j,k)}(g\circ^{(i)}f)&=\mathrm{tg}^{(j,k)}(g)\circ^{(i)}\mathrm{tg}^{(j,k)}(f).
\end{align*}
\item[(MS4)]\label{MS4} For every $0\leq k<l<m$, if $f\in C_{k}$, then
\allowdisplaybreaks
\begin{align*}
\mathrm{i}^{(m,l)}(\mathrm{i}^{(l,k)}(f))=\mathrm{i}^{(m,k)}(f).
\end{align*}
\item[(MS5)]\label{MS5} For every $0\leq k<l$, if $f\in C_{k}$, then
$$
\mathrm{sc}^{(k,l)}(\mathrm{i}^{(l,k)}(f))=f=
\mathrm{tg}^{(k,l)}(\mathrm{i}^{(l,k)}(f)).
$$
\item[(MS6)]\label{MS6} For every $0\leq j< k$, if $f\in C_{k}$, then
$$
f\circ^{(j)}
(\mathrm{i}^{(k,j)}(\mathrm{sc}^{(j,k)}(f)))
=
f
=
(\mathrm{i}^{(k,j)}(\mathrm{tg}^{(j,k)}(f)))
\circ^{(j)}
f.
$$
\item[(MS7)]\label{MS7} For every $0\leq j <k$, if $h,g,f\in C_{k}$ are such that 
\allowdisplaybreaks
\begin{align*}
\mathrm{sc}^{(j,k)}(h)&=\mathrm{tg}^{(j,k)}(g);
&
\mathrm{sc}^{(j,k)}(g)&=\mathrm{tg}^{(j,k)}(f),
\end{align*}
then
$$
h\circ^{(j)}(g\circ^{(j)}f)=
(h\circ^{(j)}g)\circ^{(j)}f.
$$
\item[(MS8)]\label{MS8} For every $0\leq j< k<l$, if $g,f\in C_{k}$ are such that $\mathrm{sc}^{(j,k)}(g)=\mathrm{tg}^{(j,k)}(f)$, then
$$
\mathrm{i}^{(l,k)}(g\circ^{(j)}f)=\mathrm{i}^{(l,k)}(g)\circ^{(j)}\mathrm{i}^{(l,k)}(f).
$$
\item[(MS9)]\label{MS9} For every $0\leq i<j<k$, if $g',g,f',f\in C_{k}$ are such that
\allowdisplaybreaks
\begin{align*}
\mathrm{sc}^{(j,k)}(g')&=\mathrm{tg}^{(j,k)}(g);
&
\mathrm{sc}^{(j,k)}(f')&=\mathrm{tg}^{(j,k)}(f);\\
\mathrm{sc}^{(i,k)}(g')&=\mathrm{tg}^{(i,k)}(f');
&
\mathrm{sc}^{(i,k)}(g)&=\mathrm{tg}^{(i,k)}(f).
\end{align*}

then,
$$
(g'\circ^{(i)}f')\circ^{(j)}(g\circ^{(i)} f)
=
(g'\circ^{(j)}g)\circ^{(i)}(f'\circ^{(j)}f)
.
$$
\end{itemize}
From now on, we will call a $\coprod_{k\in \omega}k$-indexed family
$$
\zeta = (\zeta_{j,k})_{(j,k)\in \coprod_{k\in \omega}k},
$$
like the previous one, a structure of $\omega$-category on the $\omega$-sorted set $(C_{k})_{k\in\omega}$. 

From items~(MS1) and~(MS4), it follows that the different source, target and identity mappings are completely determined by the following series of mappings.
\begin{center}
\begin{tikzpicture}
[ACliment/.style={-{To [angle'=45, length=5.75pt, width=4pt, round]}}
, scale=1, 
AClimentD/.style={dogble eqgal sign distance,
-implies
}
]

\node[] (0R) at (.2,0) [color=white] {$A^{A}_{A}$};
\node[] (1L) at (2.8,0)  [color=white] {$A^{A}_{A}$};
\node[] (1R) at (3.2,0) [color=white] {$A^{A}_{A}$};
\node[] (2L) at (5.8,0)  [color=white] {$A^{A}_{A}$};
\node[] (2R) at (6.2,0) [color=white] {$A^{A}_{A}$};
\node[] (3L) at (8.8,0)  [color=white] {$A^{A}_{A}$};
\node[] (3R) at (9.2,0) [color=white] {$A^{A}_{A}$};
\node[] (4L) at (11.8,0)  [color=white] {$A^{A}_{A}$};

\node[] (5R) at (3.2,-2) [color=white] {$A^{A}_{A}$};
\node[] (6L) at (4.3,-2)  [color=white] {$A^{A}_{A}$};
\node[] (6R) at (4.7,-2) [color=white] {$A^{A}_{A}$};
\node[] (7L) at (7.3,-2)  [color=white] {$A^{A}_{A}$};
\node[] (7R) at (7.7,-2) [color=white] {$A^{A}_{A}$};
\node[] (8L) at (10.3,-2)  [color=white] {$A^{A}_{A}$};
\node[] (8R) at (10.7,-2) [color=white] {$A^{A}_{A}$};
\node[] (9L) at (11.8,-2)  [color=white] {$A^{A}_{A}$};

\node[] (C0) at (0,0) [] {$C_{0}$};
\node[] (C1) at (3,0) [] {$C_{1}$};
\node[] (C2) at (6,0) [] {$C_{2}$};
\node[] (C3) at (9,0) [] {$C_{3}$};
\node[] () at (12,0) [] {$\cdots$};
\node[] () at (3,-2) [] {$\cdots$};
\node[] (Cpn) at (4.5,-2) [] {$C_{k-1}$};
\node[] (Cn) at (7.5,-2) [] {$C_{k}$};
\node[] (Csn) at (10.5,-2) [] {$C_{k+1}$};
\node[] () at (12,-2) [] {$\cdots$};

\draw[ACliment] (0R) to node [midway, fill=white] {$\scriptstyle \mathrm{i}^{(1,0)}$} (1L);
\draw[ACliment] (1L.north west) to node [above] {$\scriptstyle \mathrm{sc}^{(0,1)}$} (0R.north east);
\draw[ACliment] (1L.south west) to node [below] {$\scriptstyle \mathrm{tg}^{(0,1)}$} (0R.south east);

\draw[ACliment] (1R) to node [midway, fill=white] {$\scriptstyle \mathrm{i}^{(2,1)}$} (2L);
\draw[ACliment] (2L.north west) to node [above] {$\scriptstyle \mathrm{sc}^{(1,2)}$} (1R.north east);
\draw[ACliment] (2L.south west) to node [below] {$\scriptstyle \mathrm{tg}^{(1,2)}$} (1R.south east);

\draw[ACliment] (2R) to node [midway, fill=white] {$\scriptstyle \mathrm{i}^{(3,2)}$} (3L);
\draw[ACliment] (3L.north west) to node [above] {$\scriptstyle \mathrm{sc}^{(2,3)}$} (2R.north east);
\draw[ACliment] (3L.south west) to node [below] {$\scriptstyle \mathrm{tg}^{(2,3)}$} (2R.south east);

\draw[ACliment] (3R) to node [midway, fill=white] {$\scriptstyle \mathrm{i}^{(4,3)}$} (4L);
\draw[ACliment] (4L.north west) to node [above] {$\scriptstyle \mathrm{sc}^{(3,4)}$} (3R.north east);
\draw[ACliment] (4L.south west) to node [below] {$\scriptstyle \mathrm{tg}^{(3,4)}$} (3R.south east);

\draw[ACliment] (6R) to node [midway, fill=white] {$\scriptstyle \mathrm{i}^{(k,k-1)}$} (7L);
\draw[ACliment] (7L.north west) to node [above] {$\scriptstyle \mathrm{sc}^{(k-1,k)}$} (6R.north east);
\draw[ACliment] (7L.south west) to node [below] {$\scriptstyle \mathrm{tg}^{(k-1,k)}$} (6R.south east);

\draw[ACliment] (7R) to node [midway, fill=white] {$\scriptstyle \mathrm{i}^{(k+1,k)}$} (8L);
\draw[ACliment] (8L.north west) to node [above] {$\scriptstyle \mathrm{sc}^{(k,k+1)}$} (7R.north east);
\draw[ACliment] (8L.south west) to node [below] {$\scriptstyle \mathrm{tg}^{(k,k+1)}$} (7R.south east);
\end{tikzpicture}
\end{center}

%
%
%

Let us note that, for every $i$, $j$, $k\in\omega$ with $i<j<k$, for simplicity of notation, we have used the same symbol, $\circ^{(i)}$, to denote the $i$-composition in $C_{j}$ and in $C_{k}$. If necessary, this can be made more precise by writing, for $j<k$, e.g., $\circ^{k,k\times_{j}k}$ or $\circ^{\binom{k}{j}}$ instead of $\circ^{(j)}$.

Moreover, from item~(MS5), it follows that, for every $0\leq j<k$, the mapping $\mathrm{i}^{(k,j)}$ is injective.

For many-sorted $\omega$-categories $\mathsf{C}$ and $\mathsf{C}'$, a \emph{many-sorted $\omega$-functor} $F$ from $\mathsf{C}$ to $\mathsf{C}'$, denoted by $F\colon \mathsf{C}\mor\mathsf{C}'$, is an $\omega$-sorted mapping $F=(F_{k})_{k\in\omega}$ from $(C_{k})_{k\in\omega}$ to $(C'_{k})_{k\in\omega}$ 
that satisfies the following conditions.
\begin{itemize}
\item[(MSi)]\label{MSi} For every $k\in \omega$, every $0\leq j<k$, and every $f\in C_{k}$,
\allowdisplaybreaks
\begin{align*}
F_{j}(\mathrm{sc}^{(j,k)}(f))
&=
\mathrm{sc}'^{(j,k)}(F_{k}(f));
&
F_{j}(\mathrm{tg}^{(j,k)}(f))
&=
\mathrm{tg}'^{(j,k)}(F_{k}(f)).
\end{align*}
\item[(MSii)]\label{MSii} For every $k\in \omega$, every $0\leq j< k$, and every $f\in C_{j}$
\allowdisplaybreaks
\begin{align*}
F_{k}(\mathrm{i}^{(k,j)}(f))
&=
\mathrm{i}'^{(k,j)}(F_{j}(f)).
\end{align*}
\item[(MSiii)]\label{MSiii} For every $k\in \omega$ and every $0\leq j< k$, if $f,g\in C_{k}$ are such that $\mathrm{sc}^{(j,k)}(g)=\mathrm{tg}^{(j,k)}(f)$, then
$$
F_{k}(g\circ^{(j)}f)=F_{k}(g)\circ'^{(j)}F_{k}(f).
$$
\end{itemize}

Items~(MSi) and~(MSii) in the previous definition state that, in the following diagram, every square on, respectively, sources, targets and identities, commutes.
\begin{center}
\begin{tikzpicture}
[ACliment/.style={-{To [angle'=45, length=5.75pt, width=4pt, round]}}
, scale=1, 
AClimentD/.style={dogble eqgal sign distance,
-implies
}
]

\node[] (0R) at (.2,0) [color=white] {$A^{A}_{A}$};
\node[] (1L) at (2.8,0)  [color=white] {$A^{A}_{A}$};
\node[] (1R) at (3.2,0) [color=white] {$A^{A}_{A}$};
\node[] (2L) at (5.8,0)  [color=white] {$A^{A}_{A}$};
\node[] (2R) at (6.2,0) [color=white] {$A^{A}_{A}$};
\node[] (3L) at (8.8,0)  [color=white] {$A^{A}_{A}$};

\node[] (0Rp) at (.2,-2) [color=white] {$A^{A}_{A}$};
\node[] (1Lp) at (2.8,-2)  [color=white] {$A^{A}_{A}$};
\node[] (1Rp) at (3.2,-2) [color=white] {$A^{A}_{A}$};
\node[] (2Lp) at (5.8,-2)  [color=white] {$A^{A}_{A}$};
\node[] (2Rp) at (6.2,-2) [color=white] {$A^{A}_{A}$};
\node[] (3Lp) at (8.8,-2)  [color=white] {$A^{A}_{A}$};

\node[] (C0) at (0,0) [] {$C_{0}$};
\node[] (C1) at (3,0) [] {$C_{1}$};
\node[] (C2) at (6,0) [] {$C_{2}$};
\node[] (C3) at (9,0) [] {$\cdots$};
\node[] (C0p) at (0,-2) [] {$C'_{0}$};
\node[] (C1p) at (3,-2) [] {$C'_{1}$};
\node[] (C2p) at (6,-2) [] {$C'_{2}$};
\node[] (C3p) at (9,-2) [] {$\cdots$};

\draw[ACliment] (0R) to node [midway, fill=white] {$\scriptstyle \mathrm{i}^{(1,0)}$} (1L);
\draw[ACliment] (1L.north west) to node [above] {$\scriptstyle \mathrm{sc}^{(0,1)}$} (0R.north east);
\draw[ACliment] (1L.south west) to node [below] {$\scriptstyle \mathrm{tg}^{(0,1)}$} (0R.south east);

\draw[ACliment] (1R) to node [midway, fill=white] {$\scriptstyle \mathrm{i}^{(2,1)}$} (2L);
\draw[ACliment] (2L.north west) to node [above] {$\scriptstyle \mathrm{sc}^{(1,2)}$} (1R.north east);
\draw[ACliment] (2L.south west) to node [below] {$\scriptstyle \mathrm{tg}^{(1,2)}$} (1R.south east);

\draw[ACliment] (2R) to node [midway, fill=white] {$\scriptstyle \mathrm{i}^{(3,2)}$} (3L);
\draw[ACliment] (3L.north west) to node [above] {$\scriptstyle \mathrm{sc}^{(2,3)}$} (2R.north east);
\draw[ACliment] (3L.south west) to node [below] {$\scriptstyle \mathrm{tg}^{(2,3)}$} (2R.south east);

\draw[ACliment] (0Rp) to node [midway, fill=white] {$\scriptstyle \mathrm{i}'^{(1,0)}$} (1Lp);
\draw[ACliment] (1Lp.north west) to node [above] {$\scriptstyle \mathrm{sc}'^{(0,1)}$} (0Rp.north east);
\draw[ACliment] (1Lp.south west) to node [below] {$\scriptstyle \mathrm{tg}'^{(0,1)}$} (0Rp.south east);

\draw[ACliment] (1Rp) to node [midway, fill=white] {$\scriptstyle \mathrm{i}'^{(2,1)}$} (2Lp);
\draw[ACliment] (2Lp.north west) to node [above] {$\scriptstyle \mathrm{sc}'^{(1,2)}$} (1Rp.north east);
\draw[ACliment] (2Lp.south west) to node [below] {$\scriptstyle \mathrm{tg}'^{(1,2)}$} (1Rp.south east);

\draw[ACliment] (2Rp) to node [midway, fill=white] {$\scriptstyle \mathrm{i}'^{(3,2)}$} (3Lp);
\draw[ACliment] (3Lp.north west) to node [above] {$\scriptstyle \mathrm{sc}'^{(2,3)}$} (2Rp.north east);
\draw[ACliment] (3Lp.south west) to node [below] {$\scriptstyle \mathrm{tg}'^{(2,3)}$} (2Rp.south east);

\draw[ACliment, shorten <=0.2cm, shorten >=0.2cm] (C0) to node [left] {$\scriptstyle F_{0}$} (C0p);
\draw[ACliment, shorten <=0.2cm, shorten >=0.2cm] (C1) to node [right] {$\scriptstyle F_{1}$} (C1p);
\draw[ACliment, shorten <=0.2cm, shorten >=0.2cm] (C2) to node [right] {$\scriptstyle F_{2}$} (C2p);
\end{tikzpicture}
\end{center}

Many-sorted $\omega$-functors can be composed. In fact, if $F$ is a many-sorted $\omega$-functor from $\mathsf{C}$ to $\mathsf{D}$ and $G$ is a many-sorted $\omega$-functor from $\mathsf{D}$ to $\mathsf{E}$, then $G\circ F$ is a many-sorted $\omega$-functor from $\mathsf{C}$ to $\mathsf{E}$, which we call the \emph{composite} of $F$ and $G$. This composition is associative. Moreover, for every many-sorted $\omega$-category $\mathsf{C}$, there is an \emph{identity} many-sorted $\omega$-functor $(\mathsf{C},\mathrm{Id}_{C},\mathsf{C})$, denoted by $\mathrm{Id}_{\mathsf{C}}\colon \mathsf{C}\mor\mathsf{C}$, where $\mathrm{Id}_{C}$ is $(\mathrm{Id}_{C_{k}})_{k\in\omega}$, the $\omega$-sorted identity mapping on $(C_{k})_{k\in\omega}$, which is the neutral element for the composition. We denote by $\mathsf{\omega Cat}^{\mathrm{MS}}$ the category consisting of many-sorted $\omega$-categories and many-sorted $\omega$-functors.

For $n\in\omega$, a \emph{many-sorted $n$-category} is an ordered pair $\mathsf{C} = ((C_{k})_{k\in n+1},\zeta)$ consisting of an $n+1$-sorted set $(C_{k})_{k\in n+1}$ and a structure of many-sorted $n$-category 
$$
\zeta = (\zeta_{j,k})_{(j,k)\in \coprod_{k\in n+1}k}
$$ 
on $(C_{k})_{k\in n+1}$, where, for every $k\in n+1$ and every $j\in k$, $\zeta_{j,k} = (\circ^{(j)},\mathrm{sc}^{(j,k)},\mathrm{tg}^{(j,k)},\mathrm{i}^{(k,j)})$ with
\begin{enumerate}
\item $\mathrm{sc}^{(j,k)}$ and $\mathrm{tg}^{(j,k)}$ mappings from $C_{k}$ to $C_{j}$  
which send $f\in C_{k}$ to its \emph{$(j,k)$-source} and \emph{$(j,k)$-target}, respectively;

\item $\mathrm{i}^{(k,j)}$ a mapping from $C_{j}$ to $C_{k}$ 
which sends $f\in C_{j}$ to its \emph{$(k,j)$-identity}; and
%
\item $\circ^{(j)}$ a partial binary operation 
$$
\circ^{(j)}\colon C_{k}\times C_{k} \dmor C_{k},
$$ 
whose domain of definition is
$$
\mathrm{Dom}\left(\circ^{(j)}\right) = 
\left\lbrace
(f,g)\in C_{k}\times C_{k}
\,
\middle|
\,
 \mathrm{sc}^{(j,k)}(g)=\mathrm{tg}^{(j,k)}(f)
 \right\rbrace,
$$ 
i.e., the fibered product of $C_{k}$ and $C_{k}$ over $\mathrm{sc}^{(j,k)}$ and $\mathrm{tg}^{(j,k)}$, and
which sends a pair $(f,g)\in \mathrm{Dom}(\circ^{(j)})$ to its \emph{$j$-composition in $C_{k}$}
\end{enumerate}
These data are required to satisfy the same conditions as an $\omega$-category, except for the fact that all of them are truncated to the $n$-th level. 

\begin{center}
\begin{tikzpicture}
[ACliment/.style={-{To [angle'=45, length=5.75pt, width=4pt, round]}}
, scale=1, 
AClimentD/.style={dogble eqgal sign distance,
-implies
}
]

\node[] (0R) at (.2,0) [color=white] {$A^{A}_{A}$};
\node[] (1L) at (2.8,0)  [color=white] {$A^{A}_{A}$};
\node[] (1R) at (3.2,0) [color=white] {$A^{A}_{A}$};
\node[] (2L) at (5.8,0)  [color=white] {$A^{A}_{A}$};
\node[] (2R) at (6.2,0) [color=white] {$A^{A}_{A}$};
\node[] (3L) at (8.8,0)  [color=white] {$A^{A}_{A}$};

\node[] (C0) at (0,0) [] {$C_{0}$};
\node[] (C1) at (3,0) [] {$C_{1}$};
\node[] (C2) at (6,0) [] {$\cdots$};
\node[] (C3) at (9,0) [] {$C_{n}$};

\draw[ACliment] (0R) to node [midway, fill=white] {$\scriptstyle \mathrm{i}^{(1,0)}$} (1L);
\draw[ACliment] (1L.north west) to node [above] {$\scriptstyle \mathrm{sc}^{(0,1)}$} (0R.north east);
\draw[ACliment] (1L.south west) to node [below] {$\scriptstyle \mathrm{tg}^{(0,1)}$} (0R.south east);

\draw[ACliment] (1R) to node [midway, fill=white] {$\scriptstyle \mathrm{i}^{(2,1)}$} (2L);
\draw[ACliment] (2L.north west) to node [above] {$\scriptstyle \mathrm{sc}^{(1,2)}$} (1R.north east);
\draw[ACliment] (2L.south west) to node [below] {$\scriptstyle \mathrm{tg}^{(1,2)}$} (1R.south east);

\draw[ACliment] (2R) to node [midway, fill=white] {$\scriptstyle \mathrm{i}^{(n,n-1)}$} (3L);
\draw[ACliment] (3L.north west) to node [above] {$\scriptstyle \mathrm{sc}^{(n-1,n)}$} (2R.north east);
\draw[ACliment] (3L.south west) to node [below] {$\scriptstyle \mathrm{tg}^{(n-1,n)}$} (2R.south east);
\end{tikzpicture}
\end{center}

The notion of many-sorted $n$-functor between $n$-categories is entirely analogous to that of many-sorted $\omega$-functor between $\omega$-categories, but truncated to the $n$-th level. We denote by $\mathsf{nCat}^{\mathrm{MS}}$ the category consisting of many-sorted $n$-categories and many-sorted $n$-functors.
\end{definition}

\section{
\texorpdfstring
{\(\w\mathsf{Cat}^{\mathrm{MS}}\) as the projective limit of \(((\bigcdot)\mathsf{Cat}^{\mathrm{MS}},U_{\mathrm{MS}})\)}
{Many-sorted Omega-Cat as a projective limit}
}

In this section, we define, for \(m,n\in\omega\) with \(n\leq m\), the many-sorted underlying functors \(U^{(n,m)}_{\mathrm{MS}}\) and \(U^{(n,\omega)}_{\mathrm{MS}}\). Moreover, we prove that the category \(\w\mathsf{Cat}^{\mathrm{MS}}\) of many-sorted \(\omega\)-categories with many-sorted \(\omega\)-functors is the projective limit of the projective system \(((\bigcdot)\mathsf{Cat}^{\mathrm{MS}},U_{\mathrm{MS}})\) where the components of the projective system will be introduced in Remark~\ref{RUmnMSProjSys}. 

We begin by defining the notions of underlying many-sorted \(m\)-category and many-sorted underlying functors.

\begin{definition}\label{DUnderlyingMS}
Let \(\mathsf{C}=((C_{k})_{k\in n+1},\zeta)\) be a many-sorted \(n\)-category with \(\zeta=(\zeta_{j,k})_{(j,k)\in\coprod_{k\in n+1}k}\) and, for every \((j,k)\in\coprod_{k\in n+1}k\),
\[
\zeta_{j,k}=(\circ^{(j)},\mathrm{sc}^{(j,k)},\mathrm{tg}^{(j,k)},\mathrm{i}^{(k,j)}).
\]
For \(m<n\) we will call the ordered pair 
\[
((C_{k})_{k\in m+1},\zeta^{<m}),
\]
with \(\zeta^{<m}=(\zeta_{j,k})_{(j,k)\in\coprod_{k\in m+1}k}\), the \emph{underlying many-sorted \(m\)-category} of \(\mathsf{C}\) and we will denote it by \(\mathsf{C}^{<m}\). Thus, it follows that \(\mathsf{C}^{<m}\) is a many-sorted \(m\)-category.

Moreover, for a many-sorted \(n\)-functor \(F=(F_{k})_{k\in n+1}\) from the many-sorted \(n\)-category \(\mathsf{C}\) to the many-sorted \(n\)-category \(\mathsf{C}'\), we will call the family \((F_{k})_{k\in m+1}\) the \emph{underlying many-sorted \(m\)-functor} of \(F\) and we will denote it by \(F^{<m}\). Thus, it follows that \(F^{<m}\) is a many-sorted \(m\)-functor from \(\mathsf{C}^{<m}\) to \(\mathsf{C'}^{<m}\).

A similar construction applies if \(\mathsf{C}\) is a many-sorted \(\omega\)-category and \(F\) is a many-sorted \(\omega\)-functor.
\end{definition}

\begin{definition}\label{DUmnMS}
Let \(m<n\). We let \(U^{(m,n)}_{\mathrm{MS}}\) stand for the covariant functor from \(\mathsf{nCat}^{\mathrm{MS}}\) to \(\mathsf{mCat}^{\mathrm{MS}}\) defined as follows:
\begin{enumerate}
\item
its object mapping sends a many-sorted \(n\)-category \(\mathsf{C}\) to its underlying many-sorted \(m\)-category \(\mathsf{C}^{<m}\); and
\item
its morphism mapping sends a many-sorted \(n\)-functor \(F\) from \(\mathsf{C}\) to \(\mathsf{C}'\) to its underlying many-sorted \(m\)-functor \(F^{<m}\) from \(\mathsf{C}^{<m}\) to \(\mathsf{C}'^{<m}\).
\end{enumerate}

Following a similar construction, for \(n\in\omega\) we define the covariant functor \(U^{(n,\omega)}_{\mathrm{MS}}\) from \(\w\mathsf{Cat}^{\mathrm{MS}}\) to \(\mathsf{nCat}^{\mathrm{MS}}\).

Moreover, we let \(U^{(n,n)}_{\mathrm{MS}}\) stand for \(\mathrm{Id}_{\mathsf{nCat}^{\mathrm{MS}}}\) the identity functor on \(\mathsf{nCat}^{\mathrm{MS}}\).
\end{definition}

\begin{remark}\label{RUmnMSProjSys}
We let \((\bigcdot)\mathsf{Cat}^{\mathrm{MS}}\) stand for the family of categories \((\mathsf{nCat}^{\mathrm{MS}})_{n\in\omega}\) and we let \(U_{\mathrm{MS}}\) stand for the family of functors \((U^{(m,n)}_{\mathrm{MS}})_{(m,n)\in\leq}\). Following Definition~\ref{DUnderlyingMS}, for every \(m<n<p\), every many-sorted \(p\)-category \(\mathsf{C}\) and every many-sorted \(p\)-functor \(F\), it follows that \((\mathsf{C}^{<n})^{<m}=\mathsf{C}^{<m}\) and \((F^{<n})^{<m}=F^{<m}\), that is, 
\[
U^{(m,n)}_{\mathrm{MS}}\circ U^{(n,p)}_{\mathrm{MS}}=U^{(m,p)}_{\mathrm{MS}}.
\]
Therefore, we obtain the projective system \(((\bigcdot)\mathsf{Cat}^{\mathrm{MS}},U_{\mathrm{MS}})\) in \(\mathsf{Cat}\) and its projective limit \(\varprojlim ((\bigcdot)\mathsf{Cat}^{\mathrm{MS}},U_{\mathrm{MS}})\).
\end{remark}

Our next goal is to prove that such limit is isomorphic to \((\w\mathsf{Cat}^{\mathrm{MS}},U^{(\bigcdot,\omega)}_{\mathrm{MS}})\) with \(U^{(\bigcdot,\omega)}_{\mathrm{MS}}\) being the family of functors \((U^{(n,\omega)}_{\mathrm{MS}})_{n\in\omega}\). That the composition of functors \(U^{(m,n)}_{\mathrm{MS}}\circ U^{(n,\omega)}_{\mathrm{MS}}\) is equal to the functor \(U^{(m,\omega)}_{\mathrm{MS}}\) follows from the just stated consideration. Thus, all that remains to be proven is that \((\w\mathsf{Cat}^{\mathrm{MS}},U^{(\bigcdot,\omega)}_{\mathrm{MS}})\) satisfies the universal property.

\begin{proposition}
\label{PUnivPropMS}
Let \(\mathsf{L}\) be a category and, for every \(n\in\omega\), let \(F_{n}\) be a functor from \(\mathsf{L}\) to \(\mathsf{nCat}^{\mathrm{MS}}\). If, for every \(m,n\in\omega\) with \(m\leq n\), \(F_{m} = U^{(m,n)}_{\mathrm{MS}} \circ F_{n}\), then there exists a unique functor \(\langle F_{k}\rangle_{k\in\omega}\) from \(\mathsf{L}\) to the category \(\w\mathsf{Cat}^{\mathrm{MS}}\) such that, for every \(n\in\omega\), \(F_{n} = U^{(n,\omega)}_{\mathrm{MS}} \circ \langle F_{k}\rangle_{k\in\omega}\). The diagram in Figure~\ref{FUnivPropMS} depicts this situation.
\begin{figure}
\[
\xymatrix{
& \mathsf{L}\ar@{-->}[d]^{\exists!\langle F_{k}\rangle_{k\in\omega}}\ar@/_30pt/[ddl]_{F_{m}}\ar@/^30pt/[ddr]^{F_{n}} &
\\
& \w\mathsf{Cat}^{\mathrm{MS}}\ar[dr]^{U^{(n,\omega)}_{\mathrm{MS}}}\ar[dl]_{U^{(m,\omega)}_{\mathrm{MS}}} &
\\
\mathsf{mCat}^{\mathrm{MS}} && \mathsf{nCat}^{\mathrm{MS}}\ar[ll]^{U^{(m,n)}_{\mathrm{MS}}}
}
\]
\caption{Universal Property of \(\w\mathsf{Cat}^{\mathrm{MS}}\).}
\label{FUnivPropMS}
\end{figure}
\end{proposition}

\begin{proof}
Before proceeding any further let us fix some notation. For every object \(A\) in \(\mathsf{L}\), we will denote its image under the functor \(F_{n}\) by \(\mathsf{A}_{n}=((A_{n,k})_{k\in n+1},\zeta^{\mathsf{A}_{n}})\) with \(\zeta^{\mathsf{A}_{n}}=(\zeta^{\mathsf{A}_{n}}_{j,k})_{(j,k)\in\coprod_{k\in n+1}k}\) where, for every \((j,k)\in\coprod_{k\in n+1}k\),
\[
\zeta^{\mathsf{A}_{n}}_{j,k}
=
(\circ^{(j)\mathsf{A}_{n}},\mathrm{sc}^{(j,k)\mathsf{A}_{n}},\mathrm{tg}^{(j,k)\mathsf{A}_{n}},\mathrm{i}^{(k,j)\mathsf{A}_{n}}).
\]
Let us remark that \(\mathsf{A}_{n}\) is a many-sorted \(n\)-category. Moreover, for every morphism \(\varphi\) from \(A\) to \(B\) in \(\mathsf{L}\), we will denote its image under the functor \(F_{n}\) simply by \(F_{n}(\varphi)\). Let us remark that \(F_{n}(\varphi)\) is a many-sorted \(n\)-functor from \(\mathsf{A}_{n}\) to \(\mathsf{B}_{n}\), thus a family \(F_{n}(\varphi)=(\varphi_{n,k})_{k\in n+1}\) where, for every \(k\in n+1\), \(\varphi_{n,k}\) is a mapping from \(A_{n,k}\) to \(B_{n,k}\). The following diagram depicts the notation introduced so far.
\[
\xymatrixcolsep={18ex}
\xymatrixrowsep={12ex}
\scalebox{0.75}{\xymatrix{
A_{n,0}
  \ar[d]_{\varphi_{n,0}}
  \ar[r]|-{\mathrm{i}^{(1,0)\mathsf{A}_{n}}}
  &
A_{n,1}
  \ar[d]^{\varphi_{n,1}}
  \ar@<-2ex>[l]_{\mathrm{sc}^{(0,1)\mathsf{A}_{n}}}
  \ar[r]|-{\mathrm{i}^{(2,1)\mathsf{A}_{n}}}
  \ar@<2ex>[l]^{\mathrm{tg}^{(0,1)\mathsf{A}_{n}}}
  &
A_{n,2}
  \ar[d]^{\varphi_{n,2}}
  \ar@<-2ex>[l]_{\mathrm{sc}^{(1,2)\mathsf{A}_{n}}}
  \ar[r]|-{\mathrm{i}^{(3,2)\mathsf{A}_{n}}}
  \ar@<2ex>[l]^{\mathrm{tg}^{(1,2)\mathsf{A}_{n}}}
  &
\cdots
  \ar@<-2ex>[l]_{\mathrm{sc}^{(2,3)\mathsf{A}_{n}}}
  \ar[r]|-{\mathrm{i}^{(n,n-1)\mathsf{A}_{n}}}
  \ar@<2ex>[l]^{\mathrm{tg}^{(2,3)\mathsf{A}_{n}}}
  &
A_{n,n}
  \ar[d]^{\varphi_{n,n}}
  \ar@<-2ex>[l]_{\mathrm{sc}^{(n-1,n)\mathsf{A}_{n}}}
  \ar@<2ex>[l]^{\mathrm{tg}^{(n-1,n)\mathsf{A}_{n}}}
\\
B_{n,0}
  \ar[r]|-{\mathrm{i}^{(1,0)\mathsf{B}_{n}}}
  &
B_{n,1}
  \ar@<-2ex>[l]_{\mathrm{sc}^{(0,1)\mathsf{B}_{n}}}
  \ar[r]|-{\mathrm{i}^{(2,1)\mathsf{B}_{n}}}
  \ar@<2ex>[l]^{\mathrm{tg}^{(0,1)\mathsf{B}_{n}}}
  &
B_{n,2}
  \ar@<-2ex>[l]_{\mathrm{sc}^{(1,2)\mathsf{B}_{n}}}
  \ar[r]|-{\mathrm{i}^{(3,2)\mathsf{B}_{n}}}
  \ar@<2ex>[l]^{\mathrm{tg}^{(1,2)\mathsf{B}_{n}}}
  &
\cdots
  \ar@<-2ex>[l]_{\mathrm{sc}^{(2,3)\mathsf{B}_{n}}}
  \ar[r]|-{\mathrm{i}^{(n,n-1)\mathsf{B}_{n}}}
  \ar@<2ex>[l]^{\mathrm{tg}^{(2,3)\mathsf{B}_{n}}}
  &
B_{n,n}
  \ar@<-2ex>[l]_{\mathrm{sc}^{(n-1,n)\mathsf{B}_{n}}}
  \ar@<2ex>[l]^{\mathrm{tg}^{(n-1,n)\mathsf{B}_{n}}}
}}
\]

However, let us point out that, for every \(m\leq n\), since \(F_{m}=U^{(m,n)}_{\mathrm{MS}}\circ F_{n}\), it follows that, for every object \(A\) in \(\mathsf{L}\) and every morphism \(\varphi\) from \(A\) to \(B\) in \(\mathsf{L}\), \(\mathsf{A}_{m}=\mathsf{A}_{n}^{<m}\) and \(\varphi_{m}=\varphi_{n}^{<m}\). In particular, 
\begin{enumerate}
\item
for every \(k\in m+1\),\(A_{m,k}=A_{n,k}\); and
\item
for every \((j,k)\in\coprod_{k\in m+1}k\),
\begin{align*}
\circ^{(j)\mathsf{A}_{m}} &= \circ^{(j)\mathsf{A}_{n}},
&
\mathrm{sc}^{(j,k)\mathsf{A}_{m}} &= \mathrm{sc}^{(j,k)\mathsf{A}_{n}},
\\
\mathrm{tg}^{(j,k)\mathsf{A}_{m}} &= \mathrm{tg}^{(j,k)\mathsf{A}_{n}},
&
\mathrm{i}^{(k,j)\mathsf{A}_{m}} &= \mathrm{i}^{(k,j)\mathsf{A}_{n}};
\end{align*}
and
\item
for every \(k\in m+1\), \(\varphi_{m,k}=\varphi_{n,k}\).
\end{enumerate}

This allows us to drop the subscript of \(A_{n}\) and, for every object \(A\) in \(\mathsf{L}\), write its image under the functor \(F_{n}\) simply as \(\mathsf{A}_{n}=((A_{k})_{k\in n+1},\zeta^{\mathsf{A}_n})\) with \(\zeta^{\mathsf{A}_n}=(\zeta^{\mathsf{A}}_{j,k})_{(j,k)\in\coprod_{k\in n+1}k}\) where, for every \((j,k)\in\coprod_{k\in n+1}k\),
\[
\zeta^{\mathsf{A}}_{j,k}
=
(\circ^{(j)\mathsf{A}},\mathrm{sc}^{(j,k)\mathsf{A}},\mathrm{tg}^{(j,k)\mathsf{A}},\mathrm{i}^{(k,j)\mathsf{A}}).
\]
And, for every morphism \(\varphi\) from \(A\) to \(B\) in \(\mathsf{L}\), write its image under the functor \(F_{n}\) simply as \(F_{n}(\varphi)=(\varphi_{k})_{k\in n+1}\). As before, the following diagram fixes the notation.

\[
\xymatrixcolsep={18ex}
\xymatrixrowsep={12ex}
\scalebox{0.85}{\xymatrix{
A_{0}
  \ar[d]_{\varphi_{0}}
  \ar[r]|-{\mathrm{i}^{(1,0)\mathsf{A}}}
  &
A_{1}
  \ar[d]^{\varphi_{1}}
  \ar@<-2ex>[l]_{\mathrm{sc}^{(0,1)\mathsf{A}}}
  \ar[r]|-{\mathrm{i}^{(2,1)\mathsf{A}}}
  \ar@<2ex>[l]^{\mathrm{tg}^{(0,1)\mathsf{A}}}
  &
A_{2}
  \ar[d]^{\varphi_{2}}
  \ar@<-2ex>[l]_{\mathrm{sc}^{(1,2)\mathsf{A}}}
  \ar[r]|-{\mathrm{i}^{(3,2)\mathsf{A}}}
  \ar@<2ex>[l]^{\mathrm{tg}^{(1,2)\mathsf{A}}}
  &
\cdots
  \ar@<-2ex>[l]_{\mathrm{sc}^{(2,3)\mathsf{A}}}
  \ar[r]|-{\mathrm{i}^{(n,n-1)\mathsf{A}}}
  \ar@<2ex>[l]^{\mathrm{tg}^{(2,3)\mathsf{A}}}
  &
A_{n}
  \ar[d]^{\varphi_{n}}
  \ar@<-2ex>[l]_{\mathrm{sc}^{(n-1,n)\mathsf{A}}}
  \ar@<2ex>[l]^{\mathrm{tg}^{(n-1,n)\mathsf{A}}}
\\
B_{0}
  \ar[r]|-{\mathrm{i}^{(1,0)\mathsf{B}}}
  &
B_{1}
  \ar@<-2ex>[l]_{\mathrm{sc}^{(0,1)\mathsf{B}}}
  \ar[r]|-{\mathrm{i}^{(2,1)\mathsf{B}}}
  \ar@<2ex>[l]^{\mathrm{tg}^{(0,1)\mathsf{B}}}
  &
B_{2}
  \ar@<-2ex>[l]_{\mathrm{sc}^{(1,2)\mathsf{B}}}
  \ar[r]|-{\mathrm{i}^{(3,2)\mathsf{B}}}
  \ar@<2ex>[l]^{\mathrm{tg}^{(1,2)\mathsf{B}}}
  &
\cdots
  \ar@<-2ex>[l]_{\mathrm{sc}^{(2,3)\mathsf{B}}}
  \ar[r]|-{\mathrm{i}^{(n,n-1)\mathsf{B}}}
  \ar@<2ex>[l]^{\mathrm{tg}^{(2,3)\mathsf{B}}}
  &
B_{n}
  \ar@<-2ex>[l]_{\mathrm{sc}^{(n-1,n)\mathsf{B}}}
  \ar@<2ex>[l]^{\mathrm{tg}^{(n-1,n)\mathsf{B}}}
}}
\]

In order to define the functor \(\langle F_{k}\rangle_{k\in\omega}\), we will prove the following two claims.

\begin{claim}\label{CUnivPropMSA}
For every object \(A\) in \(\mathsf{L}\) we can uniquely assign a many-sorted \(\omega\)-category \(\mathsf{A}_{\omega}\) such that, for every \(n\in\omega\),
\[
\mathsf{A}_{\omega}^{<n}
=
\mathsf{A}_{n}.
\]
\end{claim}

We consider the ordered pair \(\mathsf{A}_{\omega}=((A_{k})_{k\in\omega}, \zeta^{\mathsf{A}_{\omega}})\) with \(\zeta^{\mathsf{A}_{\omega}}=(\zeta^{\mathsf{A}}_{j,k})_{(j,k)\in\coprod_{k\in\omega}k}\) where \(A_{k}\) and \(\zeta^{\mathsf{A}}_{j,k}\) are defined before.

The fact that \(\mathsf{A}_{\omega}\) is a many-sorted \(\omega\)-category follows directly from the fact that, for every \(n\in \omega\), \(\mathsf{A}_{n}\) is a many sorted \(n\)-category. Moreover, for every \(n\in\omega\), \(\mathsf{A}_{\omega}^{<n}=\mathsf{A}_{n}\). This follows from the definition of \(\mathsf{A}_{\omega}\). Therefore, all that remains to be proven is that this construction is the unique with such property.

Let \(\mathsf{C}=((C_{k})_{k\in\omega},\zeta)\) be a many-sorted \(\omega\)-category with \(\zeta=(\zeta_{j,k})_{(j,k)\in\coprod_{k\in\omega}k}\) where, for every \((j,k)\in\coprod_{k\in\omega}k\), \(\zeta_{j,k}=(\circ^{(j)},\mathrm{sc}^{(j,k)},\mathrm{tg}^{(j,k)},\mathrm{i}^{(k,j)})\) such that, for every \(n\in\omega\),
\[
\mathsf{C}^{<n}=\mathsf{A}_{n}.
\]

We want to prove that \(\mathsf{C}=\mathsf{A}_{\omega}\), that is, for every \(n\in\omega\), \(C_{n}=A_{n}\) and, for every \(n,m\in\omega\) with \(n<m\), \(\zeta_{n,m}=\zeta^{\mathsf{A}}_{n,m}\).

Let \(n\in\omega\). By assumption \(\mathsf{C}^{<n}=\mathsf{A}_{n}\), therefore \((C_{k})_{k\in n+1}=(A_{k})_{k\in n+1}\). In particular, \(C_{n}=A_{n}\).

Let \(n,m\in\omega\) with \(n<m\). By assumption \(\mathsf{C}^{<m}=\mathsf{A}_{m}\), therefore
\[
(\zeta_{j,k})_{(j,k)\in\coprod_{k\in m+1}k}=(\zeta^{\mathsf{A}}_{j,k})_{(j,k)\in\coprod_{k\in m+1}k}.
\]
Since \((n,m)\in\coprod_{k\in m+1}k\) it follows that \(\zeta_{n,m}=\zeta_{n,m}^{\mathsf{A}}\).

Claim~\ref{CUnivPropMSA} follows.

\begin{claim}\label{CUnivPropMSB}
For every morphism \(\varphi\) from \(A\) to \(B\) in \(\mathsf{L}\) we can uniquely assign a many-sorted \(\omega\)-functor \(F_{\omega}(\varphi)\) such that, for every \(n\in\omega\),
\[
F_{\omega}(\varphi)^{<n}
=
F_{n}(\varphi).
\]
\end{claim}

We consider the family \(F_{\omega}(\varphi)=(\varphi_{k})_{k\in\omega}\) where \(\varphi_{k}\) is defined before.

The fact that \(F_{\omega}(\varphi)\) is a many-sorted \(\omega\)-functor from \(\mathsf{A}_{\omega}\) to \(\mathsf{B}_{\omega}\) follows from the fact that, for every \(n\in\omega\), \(F_{n}(\varphi)\) is a many-sorted \(n\)-functor from \(\mathsf{A}_{n}\) to \(\mathsf{B}_{n}\). Moreover, for every \(n\in\omega\), \(F_{\omega}(\varphi)^{<n}=F_{n}(\varphi)\). This follows from the definition of \(F_{\omega}(\varphi)\). Therefore, all that remains to be proven is that this construction is the unique with such property.

Let \(\psi=(\psi_{k})_{k\in\omega}\) be a many-sorted \(\omega\)-functor from \(\mathsf{A}_{\omega}\) to \(\mathsf{B}_{\omega}\) such that, for every \(n\in\omega\), 
\[
\psi^{<n}
=
F_{n}(\varphi).
\]

We want to prove that \(\psi=F_{\omega}(\varphi)\), that is, for every \(n\in\omega\), \(\psi_{n}=\varphi_{n}\).

Let \(n\in\omega\). By assumption \(\psi^{<n}=F_{n}(\varphi)\), that is, \((\psi_{k})_{k\in n+1}=(\varphi_{k})_{k\in n+1}\). In particular, \(\psi_{n}=\varphi_{n}\).

Claim~\ref{CUnivPropMSB} follows. 

Therefore, we define the functor \(\langle F_{k}\rangle_{k\in\omega}\) from \(\mathsf{L}\) to \(\w\mathsf{Cat}^{\mathrm{MS}}\) whose object mapping sends an object \(A\) in \(\mathsf{L}\) to the many-sorted \(\omega\)-category \(\mathsf{A}_{\omega}\), and whose morphism mapping sends a morphism \(\varphi\) from \(A\) to \(B\) in \(\mathsf{L}\) to the many-sorted \(\omega\)-functor \(F_{\omega}(\varphi)\) from \(\mathsf{A}_{\omega}\) to \(\mathsf{B}_{\omega}\).

It is straightforward to check that \(\langle F_{k}\rangle_{k\in\omega}\) preserves compositions and identities.

Claims~\ref{CUnivPropMSA} and \ref{CUnivPropMSB} prove that the functor is well-defined and that, for every \(n\in\omega\), the diagram
\[
\xymatrix{
\mathsf{L}
  \ar[d]_{\langle F_{k}\rangle_{k\in\omega}}
  \ar[dr]^{F_{n}}
\\
\w\mathsf{Cat}^{\mathrm{MS}}
  \ar[r]_{U^{(n,\omega)}_{\mathrm{MS}}}
  &
\mathsf{nCat}^{\mathrm{MS}}
}
\]
commutes. Moreover, from the unicity of both constructions it follows the unicity of the functor \(\langle F_{k}\rangle_{k\in\omega}\).
\end{proof}

\begin{remark}
Let \(m\in\omega\). For every single-sorted \(\omega\)-category \(\mathsf{C}=(C,\xi)\), its underlying \(m\)-category presented in Definition~\ref{DUnderlying}, i.e., the image of \(\mathsf{C}\) under the object mapping of the functor \(U^{(m,\omega)}\), is denoted by \(\mathsf{C}^{<m}=(C_{m},\xi^{<m})\). Moreover, for every many-sorted \(\omega\)-category \(\mathsf{C}=((C_{k})_{k\in\omega},\zeta)\), its underlying many-sorted \(m\)-category presented in Definition~\ref{DUnderlyingMS}, i.e., the image of \(\mathsf{C}\) under the object mapping of the functor \(U^{(m,\omega)}_{\mathrm{MS}}\), is denoted by \(\mathsf{C}^{<m}=((C_{k})_{k\in m+1},\zeta^{<m})\). Similarly for \(m<n\) and the functors \(U^{(m,n)}\) and \(U^{(m,n)}_{\mathrm{MS}}\).

 Note that we are using the same notation for both definitions. Despite the abuse of notation, it will become clear from the context whether we are using the single-sorted or the many-sorted definition.
\end{remark}

\section{
\texorpdfstring
{From single-sorted to many-sorted \(n\)-categories}
{From single-sorted to many-sorted n-categories}
}
In this section, we associate to a single-sorted $n$-category a many-sorted $n$-category by means of the notion of graduation. Moreover, we also associate to a single-sorted $n$-functor between single-sorted $n$-categories a many-sorted $n$-functor between the associated many-sorted $n$-categories. We end this section by showing that the above construction is functorial.

We begin by defining  the notion of graduation of a single-sorted $n$-category.

\begin{definition}\label{DGrad} 
Let $\mathsf{C}=(C,(\xi_{k})_{k\in n})$ be a single-sorted $n$-category with $(\xi_{k})_{k\in n} = (\#^{k},\mathrm{sc}^{k}, \mathrm{tg}^{k})_{k\in n}$. For $k\in n$, we will say that $f\in C$ is a \emph{$k$-cell} if  $f=\mathrm{sc}^{k}(f)$ or, equivalently, by item~(C1) in Definition~\ref{DCat}, if $f=\mathrm{tg}^{k}(f)$. Then we denote by $C_{k}$ the set of $k$-cells of $\mathsf{C}$, i.e.,
$$
\mathrm{C}_{k}=\{f\in C\mid f=\mathrm{sc}^{k}(f)\}.
$$
Moreover, we let $C_{n}$ stand for $C$ and, by analogy to the above, the elements of $C_{n}$ will be called  \emph{$n$-cells}. 

Now, concerning the $n+1$-sorted set $(C_{k})_{k\in n+1}$, we define the following mappings.

For every $0\leq j<k\leq n$, we have that $C_{j}\subseteq C_{k}$. Indeed, let $f$ be a $j$-cell in $C_{j}$, i.e., $f=\mathrm{sc}^{j}(f)$. Then, by item~(2C1) in Definition~\ref{D2Cat}, it happens that
$$
\mathrm{sc}^{k}(f)=\mathrm{sc}^{k}(\mathrm{sc}^{j}(f))=\mathrm{sc}^{j}(f)=f.
$$
Taking this into account, we let $\mathrm{i}^{(k,j)}$ stand for $\mathrm{in}_{C_{j},C_{k}}$, the canonical inclusion of $C_{j}$ into $C_{k}$. Thus
$$
\mathrm{i}^{(k,j)}\colon C_{j}\mor C_{k}.
$$

For every $0\leq j<k\leq n$, we have that $\mathrm{sc}^{j}[C_{k}]\subseteq \mathrm{C}_{j}$. Indeed, let $f$ be a $k$-cell in $C_{k}$. Then, by item~(C1) in Definition~\ref{DCat}, it happens that
$$
\mathrm{sc}^{j}(\mathrm{sc}^{j}(f))
=\mathrm{sc}^{j}(f).
$$
Therefore, $\mathrm{sc}^{j}{\,\!\upharpoonright}_{ C_{k}}$ corestricts to $\mathrm{C}_{j}$. Then we let $\mathrm{sc}^{(j,k)}$ stand for $\mathrm{sc}^{j}\bigr|_{C_{k}}^{C_{j}}$, the corresponding birestriction. Thus
$$
\mathrm{sc}^{(j,k)}\colon C_{k}\mor C_{j}.
$$

For every $0\leq j<k\leq n$, we have that $\mathrm{tg}^{j}[C_{k}]\subseteq \mathrm{C}_{j}$. Indeed, let $f$ be  a $k$-cell in $C_{k}$. Then, by item~(C1) in Definition~\ref{DCat}, it happens that
$$ 
\mathrm{sc}^{j}(\mathrm{tg}^{j}(f))
=\mathrm{tg}^{j}(f)
$$
or, what is equivalent, that
$$
\mathrm{tg}^{j}(\mathrm{tg}^{j}(f))
=\mathrm{tg}^{j}(f).
$$
Therefore, $\mathrm{tg}^{j}{\,\!\upharpoonright}_{ C_{k}}$ corestricts to $\mathrm{C}_{j}$. Then we let $\mathrm{tg}^{(j,k)}$ stand for $\mathrm{tg}^{j}\bigr|_{C_{k}}^{C_{j}}$, the corresponding birestriction. Thus
$$
\mathrm{tg}^{(j,k)}\colon C_{k}\mor C_{j}.
$$

We summarize in the following diagram the mappings just defined.
\begin{center}
\begin{tikzpicture}
[ACliment/.style={-{To [angle'=45, length=5.75pt, width=4pt, round]}}
, scale=1, 
AClimentD/.style={dogble eqgal sign distance,
-implies
}
]

\node[] (0R) at (.2,0) [color=white] {$A^{A}_{A}$};
\node[] (1L) at (2.8,0)  [color=white] {$A^{A}_{A}$};
\node[] (1R) at (3.2,0) [color=white] {$A^{A}_{A}$};
\node[] (2L) at (5.8,0)  [color=white] {$A^{A}_{A}$};
\node[] (2R) at (6.2,0) [color=white] {$A^{A}_{A}$};
\node[] (3L) at (8.8,0)  [color=white] {$A^{A}_{A}$};

\node[] (C0) at (0,0) [] {$C_{0}$};
\node[] (C1) at (3,0) [] {$C_{1}$};
\node[] (C2) at (6,0) [] {$\cdots$};
\node[] (C3) at (9,0) [] {$C_{n}$};

\draw[ACliment] (0R) to node [midway, fill=white] {$\scriptstyle \mathrm{i}^{(1,0)}$} (1L);
\draw[ACliment] (1L.north west) to node [above] {$\scriptstyle \mathrm{sc}^{(0,1)}$} (0R.north east);
\draw[ACliment] (1L.south west) to node [below] {$\scriptstyle \mathrm{tg}^{(0,1)}$} (0R.south east);

\draw[ACliment] (1R) to node [midway, fill=white] {$\scriptstyle \mathrm{i}^{(2,1)}$} (2L);
\draw[ACliment] (2L.north west) to node [above] {$\scriptstyle \mathrm{sc}^{(1,2)}$} (1R.north east);
\draw[ACliment] (2L.south west) to node [below] {$\scriptstyle \mathrm{tg}^{(1,2)}$} (1R.south east);

\draw[ACliment] (2R) to node [midway, fill=white] {$\scriptstyle \mathrm{i}^{(n,n-1)}$} (3L);
\draw[ACliment] (3L.north west) to node [above] {$\scriptstyle \mathrm{sc}^{(n-1,n)}$} (2R.north east);
\draw[ACliment] (3L.south west) to node [below] {$\scriptstyle \mathrm{tg}^{(n-1,n)}$} (2R.south east);
\end{tikzpicture}
\end{center}


Moreover, for every $0\leq j<k\leq n$, if $f,g\in C_{k}$ are such that $\mathrm{sc}^{j}(g)=\mathrm{tg}^{j}(f)$, then we have that the $j$-composite $g\#^{j}f$ belongs to $C_{k}$. Indeed, by item~(2C2) in Definition~\ref{D2Cat} and taking into account that $f$ and $g$ are $k$-cells, it happens that
$$
\mathrm{sc}^{k}(g\#^{j}f)=\mathrm{sc}^{k}(g)\#^{j}\mathrm{sc}^{k}(f)=g\#^{j}f.
$$
Therefore, $\#^{j}{\,\!\upharpoonright}_{ C_{k}\times C_{k}}$ corestricts to $\mathrm{C}_{k}$.
Then we let $\mathrm{\circ}^{(j)}$ stand for $\#^{j}\bigr|_{ C_{k}\times C_{k}}^{C_{k}}$, the corresponding birestriction. Thus
$$
\mathrm{\circ}^{(j)}\colon C_{k}\times C_{k}\dmor C_{k}.
$$
Note that, for every $f,g\in C_{k}$ the composition $g\circ^{(j)}f$ is defined if, and only if, $g\#^{j}f$ is defined, which, in its turn, by item~(C2) in Definition~\ref{DCat}, is fulfilled if, and only if, $\mathrm{sc}^{j}(g)=\mathrm{tg}^{j}(f)$. Thus, taking into account item~(2C1) in Definition~\ref{D2Cat} and the fact that $g,f$ are $k$-cells, we have that
\allowdisplaybreaks
\begin{align*}
\mathrm{sc}^{j}(g)&=\mathrm{tg}^{j}(f)
&\mbox{if, and only if,}&&
\mathrm{sc}^{(j,k)}(g)&=\mathrm{tg}^{(j,k)}(f).
\end{align*}

We will call the ordered pair 
$$
\left(\left(
C_{k}
\right)_{k\in n+1}, 
\left(
\circ^{(j)},
\mathrm{sc}^{(j,k)},
\mathrm{tg}^{(j,k)},\mathrm{i}^{(k,j)}
\right)_{(j,k)\in \coprod_{k\in n+1}k}
\right)
$$
the \emph{graduation of $\mathsf{C}$} and we will denote it by $\mathrm{Gd}^{(n)}(\mathsf{C})$. 
\end{definition}

In the following proposition we prove that the graduation of a single-sorted $n$-category is a many-sorted $n$-category.

\begin{proposition}\label{PGrad} 
Let $\mathsf{C}=(C,(\xi_{k})_{k\in n})$ be a single-sorted $n$-category with $(\xi_{k})_{k\in n} = (\#^{k},\mathrm{sc}^{k}, \mathrm{tg}^{k})_{k\in n}$. Then $\mathrm{Gd}^{(n)}(\mathsf{C})$ is a many-sorted $n$-category.
\end{proposition}

\begin{proof}
We want to prove that $ (\circ^{(j)},\mathrm{sc}^{(j,k)},\mathrm{tg}^{(j,k)},\mathrm{i}^{(k,j)})_{(j,k)\in \coprod_{k\in n+1}k}$, as described in Definition~\ref{DGrad}, is a structure of many-sorted $n$-category on $(C_{k})_{k\in n+1}$. Thus, we must check that it satisfies the conditions stated in Definition~\ref{DnCatMS}.

\textsf{(MS1)} For every $0\leq i<j<k\leq n$, if $f\in C_{k}$, then  the following chain of equalities holds
\allowdisplaybreaks
\begin{align*}
\mathrm{sc}^{(i,j)}\left(
\mathrm{sc}^{(j,k)}\left(
f
\right)\right)
&=
\mathrm{sc}^{i}\left(
\mathrm{sc}^{j}\left(
f
\right)\right)
\tag{1}
\\&=
\mathrm{sc}^{i}\left(
\mathrm{tg}^{j}\left(
f
\right)\right)
\tag{2}
\\&=
\mathrm{sc}^{(i,j)}\left(
\mathrm{tg}^{(j,k)}\left(
f
\right)\right).
\tag{3}
\end{align*}

The first equality unravels the definition of the mappings $\mathrm{sc}^{(i,j)}$ and $\mathrm{sc}^{(j,k)}$ presented in Definition~\ref{DGrad}; the second equality follows from item~(2C1) in Definition~\ref{D2Cat}; finally, the last equality unravels  the definition of the mappings $\mathrm{sc}^{(i,j)}$ and $\mathrm{tg}^{(j,k)}$ presented in Definition~\ref{DGrad}. 

By a similar argument we also have that 
\allowdisplaybreaks
\begin{align*}
\mathrm{sc}^{(i,j)}\left(
\mathrm{sc}^{(j,k)}\left(
f
\right)\right)&=
\mathrm{sc}^{(i,k)}\left(
f
\right).
\\
\mathrm{tg}^{(i,j)}\left(
\mathrm{sc}^{(j,k)}\left(
f
\right)\right)
&=
\mathrm{tg}^{(i,k)}\left(
f
\right)=
\mathrm{tg}^{(i,j)}\left(
\mathrm{tg}^{(j,k)}\left(
f
\right)\right).
\end{align*}

\textsf{(MS2)} For every $0\leq j< k\leq n$, if $g,f\in C_{k}$ are such that $\mathrm{sc}^{(j,k)}(g)=\mathrm{tg}^{(j,k)}(f)$, then the following chain of equalities holds
\allowdisplaybreaks
\begin{align*}
\mathrm{sc}^{(j,k)}\left(
g\circ^{(j)}f
\right)&=
\mathrm{sc}^{j}\left(
g\#^{j}f
\right)
\tag{1}
\\&=
\mathrm{sc}^{j}\left(
f
\right)
\tag{2}
\\&=
\mathrm{sc}^{(j,k)}\left(
f
\right).
\tag{3}
\end{align*}

The first equality unravels the definition of the mappings $\mathrm{sc}^{(j,k)}$ and $\circ^{(j)}$ presented in Definition~\ref{DGrad}; the second equality follows from item~(C3) in Definition~\ref{DCat}; finally, the last equality unravels the definition of the mapping $\mathrm{sc}^{(j,k)}$ presented in Definition~\ref{DGrad}.

By a similar argument we also have that 
$$
\mathrm{tg}^{(j,k)}\left(
g\circ^{(j)}f
\right)
=
\mathrm{tg}^{(j,k)}\left(
g
\right).
$$

\textsf{(MS3)} For every $0\leq i<j< k\leq n$, if $g,f\in C_{k}$ are such that $\mathrm{sc}^{(i,k)}(g)=\mathrm{tg}^{(i,k)}(f)$ then the following chain of equalities holds
\allowdisplaybreaks
\begin{align*}
\mathrm{sc}^{(j,k)}\left(
g\circ^{(i)}f
\right)&=
\mathrm{sc}^{j}\left(
g\#^{i}f
\right)
\tag{1}
\\&=
\mathrm{sc}^{j}\left(
g\right)
\#^{i}\mathrm{sc}^{j}\left(
f
\right)
\tag{2}
\\&=
\mathrm{sc}^{(j,k)}\left(
g
\right)
\circ^{(i)}
\mathrm{sc}^{(j,k)}\left(
f
\right).
\tag{3}
\end{align*}

The first equality unravels the definition of the mappings $\mathrm{sc}^{(j,k)}$ and $\circ^{(i)}$ presented in Definition~\ref{DGrad}; the second equality follows from item~(2C2) in Definition~\ref{D2Cat}; finally, the last equality unravels the definition of the mappings $\mathrm{sc}^{(j,k)}$ and $\circ^{(i)}$ presented in Definition~\ref{DGrad}.

By a similar argument we also have that 
$$
\mathrm{tg}^{(j,k)}\left(
g\circ^{(i)}f
\right)
=
\mathrm{tg}^{(j,k)}\left(
g
\right)
\circ^{(i)}
\mathrm{tg}^{(j,k)}\left(
f
\right).
$$

\textsf{(MS4)}  For every $0\leq k<l<m\leq n$, if $f\in C_{k}$, then the following equality holds
\allowdisplaybreaks
\begin{align*}
\mathrm{i}^{(m,l)}\left(
\mathrm{i}^{(l,k)}\left(
f
\right)\right)=\mathrm{i}^{(m,k)}\left(
f
\right).
\end{align*}
Let us recall that all the mappings under consideration are inclusions.

\textsf{(MS5)} For every $0\leq k<l\leq n$, if $f\in C_{k}$, then the following chain of equalities holds
\allowdisplaybreaks
\begin{align*}
\mathrm{sc}^{(k,l)}\left(
\mathrm{i}^{(l,k)}(f)
\right)
&=
\mathrm{sc}^{k}(f)
\tag{1}
\\&=
f.
\tag{2}
\end{align*}

The first equality unravels the definition of the mappings $\mathrm{sc}^{(k,l)}$ and $\mathrm{i}^{(l,k)}$ presented in Definition~\ref{DGrad}; finally, for the last equality, it suffices to take into account the fact that $f\in C_{k}$, thus $f=\mathrm{sc}^{k}(f)$.

By a similar argument we also have that 
$$
\mathrm{tg}^{(k,l)}\left(
\mathrm{i}^{(l,k)}(f)
\right)
=
f.
$$

\textsf{(MS6)} For every $0\leq j<k$, if $f\in C_{k}$, the following chain of equalities holds
\allowdisplaybreaks
\begin{align*}
f\circ^{(j)}
\left(
\mathrm{i}^{(k,j)}\left(
\mathrm{sc}^{(j,k)}(f)
\right)\right)
&=
f\#^{j}
\mathrm{sc}^{j}(f)
\tag{1}
\\&=
f.
\tag{2}
\end{align*}

The first equality unravels the definition of the mappings $\circ^{(j)}$, $\mathrm{sc}^{(j,k)}$ and $\mathrm{i}^{(k,j)}$ presented in Definition~\ref{DGrad}; the last equality follows from item~(C4) in Definition~\ref{DCat}.

By a similar argument we also have that 
$$
\left(
\mathrm{i}^{(k,j)}\left(
\mathrm{tg}^{(j,k)}(f)
\right)\right)
\circ^{(j)} f
=
f.
$$

\textsf{(MS7)} For every $0\leq j<k\leq n$, if $h,g,f\in C_{k}$ are such that $\mathrm{sc}^{(j,k)}(g)=\mathrm{tg}^{(j,k)}(f)$ and $\mathrm{sc}^{(j,k)}(h)=\mathrm{tg}^{(j,k)}(g)$, then the following chain of equalities holds
\allowdisplaybreaks
\begin{align*}
h\circ^{(j)}\left(
g\circ^{(j)}f
\right)&=
h\#^{j}\left(
g\#^{j}f
\right)
\tag{1}
\\&=
\left(
h\#^{j}g
\right)\#^{j}f
\tag{2}
\\&=
(h\circ^{(j)}g)\circ^{(j)}f.
\tag{3}
\end{align*}

The first equality unravels the definition of the composition mapping $\circ^{(j)}$ presented in Definition~\ref{DGrad}; the second equality follows from item~(C5) in Definition~\ref{DCat}; finally, the last equality unravels the definition of the composition mapping $\circ^{(j)}$ presented in Definition~\ref{DGrad}.

\textsf{(MS8)} For every $0\leq j<k<l \leq n$, if $g,f\in C_{k}$ are such that $\mathrm{sc}^{(j,k)}(g)=\mathrm{tg}^{(j,k)}(f)$, then the following chain of equalities holds
\allowdisplaybreaks
\begin{align*}
\mathrm{i}^{(l,k)}\left(
g\circ^{(j)}f
\right)&=
g\#^{j}f
\tag{1}
\\&=
\mathrm{i}^{(l,k)}(g)
\circ^{(j)} \mathrm{i}^{(l,k)}(f).
\tag{2}
\end{align*}

The first equality unravels the definition of the mappings $\mathrm{i}^{(l,k)}$ and $\circ^{(j)}$ presented in Definition~\ref{DGrad}. Note that $\mathrm{i}^{(l,k)}$ is just an inclusion; finally, the last equality  unravels the definition of the mappings $\mathrm{i}^{(l,k)}$ and $\circ^{(j)}$ as stated in Definition~\ref{DGrad}. 

\textsf{(MS9)} For every $0\leq i<j<k\leq n$, if $g',g,f',f\in C_{k}$ are such that 
\allowdisplaybreaks
\begin{align*}
\mathrm{sc}^{(j,k)}(g')&=\mathrm{tg}^{(j,k)}(g);
&
\mathrm{sc}^{(j,k)}(f')&=\mathrm{tg}^{(j,k)}(f);
\\
\mathrm{sc}^{(i,k)}(g')&=\mathrm{tg}^{(i,k)}(f');
&
\mathrm{sc}^{(i,k)}(g)&=\mathrm{tg}^{(i,k)}(f),
\end{align*}
then the following chain of equalities holds
\allowdisplaybreaks
\begin{align*}
\left(g'\circ^{(i)}f'\right)\circ^{(j)}\left(
g\circ^{(i)}f\right)&=
\left(
g'\#^{i}f'
\right)\#^{j}\left(
g\#^{i}f
\right)
\tag{1}
\\&=
\left(
g'\#^{j}g
\right)\#^{i}\left(
f'\#^{j}f
\right)
\tag{2}
\\&=
\left(
g'\circ^{(j)}g
\right)\circ^{(i)}\left(
f'\circ^{(j)}f\right).
\tag{3}
\end{align*}

The first equality unravels the definition of the composition mappings $\circ^{(i)}$ and $\circ^{(j)}$ presented in Definition~\ref{DGrad}; the second equality follows from item~(2C3) in Definition~\ref{D2Cat}; finally, the last equality unravels the definition of the composition mappings $\circ^{(i)}$ and $\circ^{(j)}$ presented in Definition~\ref{DGrad}.

All in all, we can affirm that $\mathrm{Gd}^{(n)}(\mathsf{C})$ is a many-sorted $n$-category. 

This completes the proof.
\end{proof}

In the following definition we introduce the notion of a graduation of a single-sorted $n$-functor.

\begin{definition}\label{DGradFun} 
Let $\mathsf{C}=(C,(\xi_{k})_{k\in n})$, with $(\xi_{k})_{k\in n} = (\#^{k},\mathrm{sc}^{k},\mathrm{tg}^{k})_{k\in n}$, and $\mathsf{C}'=(C',(\xi'_{k})_{k\in n})$, with $(\xi'_{k})_{k\in n}=(\#'^{k},\mathrm{sc}'^{k},\mathrm{tg}'^{k})_{k\in n}$, be two single-sorted $n$-categories, $F\colon \mathsf{C}\mor \mathsf{C}'$ a single-sorted $n$-functor from $\mathsf{C}$ to $\mathsf{C}'$, and 
\allowdisplaybreaks
\begin{align*}
\mathrm{Gd}^{(n)}\left(\mathsf{C}\right)&=
\left(\left(
C_{k}
\right)_{k\in n+1}, 
\left(
\circ^{(j)},\mathrm{sc}^{(j,k)},\mathrm{tg}^{(j,k)},\mathrm{i}^{(k,j)}
\right)_{(j,k)\in \coprod_{k\in n+1}k}
\right) \text{ and }
\\
\mathrm{Gd}^{(n)}\left(\mathsf{C}'\right)&=
\left(\left(
C'_{k}
\right)_{k\in n+1},
\left(
\circ'^{(j)},\mathrm{sc}'^{(j,k)},\mathrm{tg}'^{(j,k)},\mathrm{i}'^{(k,j)}
\right)_{(j,k)\in \coprod_{k\in n+1}k}
\right),
\end{align*}
the graduations corresponding to $\mathsf{C}$ and $\mathsf{C}'$, respectively. Then, for every $0\leq k\leq n$, $F{\,\!\upharpoonright}_{ C_{k}}$ corestricts to $C'_{k}$. Indeed, let $f$ be a $k$-cell in $C_{k}$. Then the following chain of equalities holds
\allowdisplaybreaks
\begin{align*}
\mathrm{sc}'^{k}\left(F(f)\right)&=
F\left(
\mathrm{sc}^{k}(f)
\right)
\tag{1}
\\&=
F(f).
\tag{2}
\end{align*}

The first equality follows from Definition~\ref{DnCatSS}; the last equality follows from the fact that $f$ is a $k$-cell in $C_{k}$. Note that, for $k=n$, $F{\,\!\upharpoonright}_{ C_{n}}=F$, which, trivially, corestricts to $C'$.

We will call the $n+1$-sorted mapping $(F_{k})_{k\in n+1}$, where, for every $k\in n+1$, $F_{k}$ is the birestriction of $F$ to $C_{k}$ and $C'_{k}$, the \emph{graduation of $F$} and we will denote it by $\mathrm{Gd}^{(n)}(F)$. Thus, $\mathrm{Gd}^{(n)}(F)$ is an $n+1$-sorted mapping from $\mathrm{Gd}^{(n)}(\mathsf{C})$ to $\mathrm{Gd}^{(n)}(\mathsf{C}')$.
\end{definition}

In the following proposition we prove that the graduation of a single-sorted $n$-functor is a many-sorted $n$-functor.

\begin{proposition}\label{PGradFun} 
Let $\mathsf{C}=(C,(\xi_{k})_{k\in n})$, with $(\xi_{k})_{k\in n} = (\#^{k},\mathrm{sc}^{k},\mathrm{tg}^{k})_{k\in n}$, and $\mathsf{C}'=(C',(\xi'_{k})_{k\in n})$, with $(\xi'_{k})_{k\in n}=(\#'^{k},\mathrm{sc}'^{k},\mathrm{tg}'^{k})_{k\in n}$, be two single-sorted $n$-categories and $F\colon \mathsf{C}\mor \mathsf{C}'$ a single-sorted $n$-functor from $\mathsf{C}$ to $\mathsf{C}'$. Then $\mathrm{Gd}^{(n)}(F)$ is a many-sorted $n$-functor from $\mathrm{Gd}^{(n)}(\mathsf{C})$ to $\mathrm{Gd}^{(n)}(\mathsf{C}')$.
\end{proposition}

\begin{proof}
It suffices to prove that the many-sorted mapping $\mathrm{Gd}^{(n)}(F)$ described in Definition~\ref{DGradFun} satisfies the conditions of Definition~\ref{DnCatMS}, thus defining a many-sorted $n$-functor from $\mathrm{Gd}^{(n)}(\mathsf{C})$ to $\mathrm{Gd}^{(n)}(\mathsf{C}')$.

\textsf{(MSi)} For every $0\leq j< k\leq n$, if $f\in C_{k}$, the following chain of equalities holds
\allowdisplaybreaks
\begin{align*}
F_{j}\left(
\mathrm{sc}^{(j,k)}(f)
\right)
&=
F\left(
\mathrm{sc}^{j}(f)
\right)
\tag{1}
\\&=
\mathrm{sc}'^{j}\left(
F(f)
\right)
\tag{2}
\\&=
\mathrm{sc}'^{(j,k)}\left(
F_{k}(f)
\right).
\tag{3}
\end{align*}

The first equality unravels the definition of the mapping $\mathrm{sc}^{(j,k)}$ presented in Definition~\ref{DGrad}; the second equality follows from the fact that $F$ is a single-sorted $n$-functor; finally, the last equality unravels the definition of the mapping $\mathrm{sc}^{(j,k)}$ presented in Definition~\ref{DGrad}.

By a similar argument we also have that 
$$
F_{j}\left(\mathrm{tg}^{(j,k)}(f)\right)
=
\mathrm{tg}'^{(j,k)}\left(
F_{k}(f)
\right).
$$

\textsf{(MSii)} For every  $0\leq k<l\leq n$, if $f\in C_{k}$, the following chain of equalities holds
\allowdisplaybreaks
\begin{align*}
F_{l}\left(
\mathrm{i}^{(l,k)}(f)
\right)
&=
F(f)
\tag{1}
\\&=
\mathrm{i}'^{(l,k)}\left(
F_{k}(f)\right).
\tag{2}
\end{align*}

The first equality unravels the definition of the mapping $\mathrm{i}^{(l,k)}$ presented in Definition~\ref{DGrad}. Let us recall that $\mathrm{i}^{(l,k)}$ is just an inclusion mapping; finally, the last equality unravels the definition of the mapping $\mathrm{i}'^{(l,k)}$ presented in Definition~\ref{DGrad} which, again, is just an inclusion mapping.

\textsf{(MSiii)} For every $0\leq j<k\leq n$, if $f,g\in C_{k}$ are such that $\mathrm{sc}^{(j,k)}(g)=\mathrm{tg}^{(j,k)}(f)$, then the following chain of equalities holds
\allowdisplaybreaks
\begin{align*}
F_{k}\left(
g\circ^{(j)}f
\right)
&=
F\left(
g\#^{j}f
\right)
\tag{1}
\\&=
F\left(
g
\right)\#^{j}F\left(
f
\right)
\tag{2}
\\&=
F_{k}(g)
\circ^{(j)}
F_{k}(f).
\tag{3}
\end{align*}

The first equality unravels the definition of the composition mapping $\circ^{(j)}$ presented in Definition~\ref{DGrad}; the second equality follows from the fact that $F$ is a single-sorted $n$-functor; finally, the last equality  unravels the definition of the composition mapping $\circ^{(j)}$ presented in Definition~\ref{DGrad} taking into consideration the fact that $f$ and $g$ are $k$-cells.

Therefore, $\mathrm{Gd}^{(n)}(F)$ is a many-sorted $n$-functor from $\mathrm{Gd}^{(n)}(\mathsf{C})$ to $\mathrm{Gd}^{(n)}(\mathsf{C}')$.
\end{proof}

We are now in position to prove that the graduation construction is functorial.

\begin{proposition}\label{PGradFunctor} 
$\mathrm{Gd}^{(n)}\colon \mathsf{nCat}\mor \mathsf{nCat}^{\mathrm{MS}}$ is a functor.
\end{proposition}

\begin{proof}
From Propositions~\ref{PGrad} and~\ref{PGradFun}, it suffices to prove that $\mathrm{Gd}^{(n)}$ preserves single-sorted identity $n$-functors and the composition of single-sorted $n$-functors.

Let $\mathsf{C}$ be a single-sorted $n$-category and consider the single-sorted identity $n$-functor $\mathrm{Id}_{\mathsf{C}}\colon\mathsf{C}\mor\mathsf{C}$, then the following chain of equalities holds
\begin{align*}
\mathrm{Gd}^{(n)}\left(
\mathrm{Id}_{\mathsf{C}}
\right)
&=
\left(
\mathrm{Id}_{k}
\right)_{k\in n+1}
\tag{1}
\\&=
\left(
\mathrm{Id}_{C_{k}}
\right)_{k\in n+1}
\tag{2}
\\&=
\mathrm{Id}_{\mathrm{Gd}^{(n)}(\mathsf{C})}.
\tag{3}
\end{align*}

The first equality unpacks Definition~\ref{DGradFun} on the graduation of $\mathrm{Id}_{\mathsf{C}}$; the second equality follows from the fact that the identity mapping on $C$, when restricted to $C_{k}$,  retrieves the identity mapping on $C_{k}$, for every $k\in n+1$; finally, the last equality recovers the definition of the many-sorted identity $n$-functor on $\mathrm{Gd}^{(n)}(\mathsf{C})$ stated in Definition~\ref{DnCatMS}.

Thus, $\mathrm{Gd}^{(n)}$ preserves identities.

Now let $F\colon \mathsf{C}\mor \mathsf{C}'$ and $G\colon\mathsf{C}'\mor\mathsf{C}''$ be single-sorted $n$-functors between single-sorted $n$-categories. Then we want to check that 
$$
\mathrm{Gd}^{(n)}\left(
G\circ F
\right)=\mathrm{Gd}^{(n)}\left(
G
\right)\circ\mathrm{Gd}^{(n)}\left(
F
\right).
$$
For $0\leq k\leq n$, we have the following chain of equalities 
\allowdisplaybreaks
\begin{align*}
\mathrm{Gd}^{(n)}\left(
G\circ F\right)_{k}&=
\left(
G\circ F\right)_{k}
\tag{1}
\\&=
G_{k}\circ F_{k}
\tag{2}
\\&=
\mathrm{Gd}^{(n)}(G)_{k}\circ\mathrm{Gd}^{(n)}(F)_{k}
\tag{3}
\\&=
\left(
\mathrm{Gd}^{(n)}(G)\circ\mathrm{Gd}^{(n)}(F)
\right)_{k}.
\tag{4}
\end{align*}

The first equality unravels Definition~\ref{DGradFun} on the $k$-th sort of $\mathrm{Gd}^{(n)}(G\circ F)$; the second equality follows from the fact that the birestriction of the composition of two mappings is the composition of its respective birestrictions. Here we take into account that, by Proposition~\ref{DGradFun}, $F{\,\!\upharpoonright_{ C_{k}}}$ corestricts to $C'_{k}$; the third equality recovers Definition~\ref{DGradFun} on the $k$-th sort of, respectively, $\mathrm{Gd}^{(n)}(G)$ and $\mathrm{Gd}^{(n)}(F)$; 
finally, the last equality recovers the definition of $k$-th sort of the composition of the many-sorted mapping $\mathrm{Gd}^{(n)}(G)\circ\mathrm{Gd}^{(n)}(F)$;

Thus, $\mathrm{Gd}^{(n)}$ preserves compositions.

This completes the proof.
\end{proof}

\section{
\texorpdfstring
{From many-sorted to single-sorted \(n\)-categories}
{From many-sorted to single-sorted n-categories}
}
In this section, we associate to a many-sorted $n$-category a single-sorted $n$-category by means of the notion of unification. Moreover, we also associate to a many-sorted $n$-functor between many-sorted $n$-categories a single-sorted $n$-functor between the associated single-sorted $n$-categories. We end this section by showing that the above construction is functorial.

We begin by introducing  the notion of unification of a many-sorted $n$-category.

\begin{definition}\label{DUnif} 
Let $\mathsf{C} = ((C_{k})_{k\in n+1},\zeta)$ be a many-sorted $n$-category  with 
$
\zeta = (\zeta_{j,k})_{(j,k)\in \coprod_{k\in n+1}k}
$ 
and, for every $(j,k)\in \coprod_{k\in n+1}k$,
$$
\zeta_{j,k} =  \left(
\circ^{(j)},\mathrm{sc}^{(j,k)},\mathrm{tg}^{(j,k)},\mathrm{i}^{(k,j)}
\right).
$$
Then, for every $k\in n$, we let $\xi_{k}$ stand for $(\#^{k},\mathrm{sc}^{k},\mathrm{tg}^{k})$, where $\#^{k}=\circ^{(k)}$, which is a partial binary operation on $C_{n}$, $\mathrm{sc}^{k} = \mathrm{i}^{(n,k)}\circ\mathrm{sc}^{(k,n)}$, which is a unary operation on $C_{n}$, and $\mathrm{tg}^{k}=
\mathrm{i}^{(n,k)}\circ\mathrm{tg}^{(k,n)}$, which is also a unary operation on $C_{n}$. 

We will call the ordered pair
$$
\left(
C_{n}, \left(
\#^{k},\mathrm{sc}^{k},\mathrm{tg}^{k}
\right)_{k\in n}
\right)
$$
the \emph{unification of $\mathsf{C}$} and we will denote it by $\mathrm{Uf}^{(n)}(\mathsf{C})$.
\end{definition}

Before continuing with the development of this section, we state in the following corollary a relation between the sources and targets of a many-sorted $n$-category and the sources and targets of its unification.

\begin{corollary}\label{CUnifScTg} 
Let $\mathsf{C} = ((C_{k})_{k\in n+1},\zeta)$ be a many-sorted $n$-category  with 
$
\zeta = (\zeta_{j,k})_{(j,k)\in \coprod_{k\in n+1}k}
$ 
and, for every $(j,k)\in \coprod_{k\in n+1}k$,
$$
\zeta_{j,k} =  \left(
\circ^{(j)},\mathrm{sc}^{(j,k)},\mathrm{tg}^{(j,k)},\mathrm{i}^{(k,j)}
\right),
$$
and $\mathrm{Uf}^{(n)}(\mathsf{C})=(C_{n}, (\#^{k},\mathrm{sc}^{k},\mathrm{tg}^{k})_{k\in n})$ its unification. Then, for every $f,g\in C_{n}$ and every $k\in n$, 
\begin{align*}
\mathrm{sc}^{(k,n)}(g)&=\mathrm{tg}^{(k,n)}(f)
&\mbox{ if, and only if,}&&
\mathrm{sc}^{k}(g)&=\mathrm{tg}^{k}(f).
\end{align*}
\end{corollary}

\begin{proof}
Let us assume that $\mathrm{sc}^{(k,n)}(g)=\mathrm{tg}^{(k,n)}(f)$. Then we can apply $\mathrm{i}^{(n,k)}$ to  both sides to infer that
$$
\mathrm{sc}^{k}(g)=\left(
\mathrm{i}^{(n,k)}\circ\mathrm{sc}^{(k,n)}\right)(g)
=
\left(\mathrm{i}^{(n,k)}\circ\mathrm{tg}^{(k,n)}\right)(f)=\mathrm{tg}^{k}(f).
$$

Now, if $\mathrm{sc}^{k}(g)=\mathrm{tg}^{k}(f)$, then the following chain of equalities holds
\allowdisplaybreaks
\begin{align*}
\mathrm{sc}^{(k,n)}(g)
&=
\mathrm{sc}^{(k,n)}
\left(
\mathrm{i}^{(n,k)}
\left(
\mathrm{sc}^{(k,n)}
\left(
g
\right)\right)\right)
\tag{1}
\\&=
\mathrm{sc}^{(k,n)}\left(
\mathrm{sc}^{k}\left(
g
\right)\right)
\tag{2}
\\&=
\mathrm{sc}^{(k,n)}\left(\mathrm{tg}^{k}\left(f
\right)\right)
\tag{3}
\\&=
\mathrm{sc}^{(k,n)}\left(
\mathrm{i}^{(n,k)}\left(
\mathrm{tg}^{(k,n)}\left(
f
\right)\right)\right)
\tag{4}
\\&=
\mathrm{tg}^{(k,n)}(f).
\tag{5}
\end{align*}

The first equality follows from the fact that, by item~(MS5) in Definition~\ref{DnCatMS}, we have that 
$
\mathrm{sc}^{(k,n)}\circ\mathrm{i}^{(n,k)}=\mathrm{Id}_{C_{k}}
$; the second equality recovers the definition of $\mathrm{sc}^{k}$ introduced in Definition~\ref{DUnif}; since we are assuming that  $\mathrm{sc}^{k}(g)=\mathrm{tg}^{k}(f)$,  
we can apply $\mathrm{sc}^{(k,n)}$ to both sides to infer the third equality; the fourth equality unravels the definition of $\mathrm{tg}^{k}$ introduced in Definition~\ref{DUnif}; finally, the last equality follows from the fact that,  by item~(MS5) in Definition~\ref{DnCatMS}, we have that 
$
\mathrm{sc}^{(k,n)}\circ\mathrm{i}^{(n,k)}=\mathrm{Id}_{C_{k}}.
$

This completes the proof.
\end{proof}

In the following proposition we prove that the unification of a many-sorted $n$-category is a single-sorted $n$-category.

\begin{proposition}\label{PUnif} 
Let $\mathsf{C} = ((C_{k})_{k\in n+1},\zeta)$ be a many-sorted $n$-category  with 
$
\zeta = (\zeta_{j,k})_{(j,k)\in \coprod_{k\in n+1}k}
$ 
and, for every $(j,k)\in \coprod_{k\in n+1}k$,
$$
\zeta_{j,k} =  \left(
\circ^{(j)},\mathrm{sc}^{(j,k)},\mathrm{tg}^{(j,k)},\mathrm{i}^{(k,j)}
\right).
$$
Then $\mathrm{Uf}^{(n)}(\mathsf{C})$ is a single-sorted $n$-category.
\end{proposition}

\begin{proof}
We want to prove that $(\#^{k},\mathrm{sc}^{k},\mathrm{tg}^{k})_{k\in n}$ as described in Definition~\ref{DUnif}, is a structure of single-sorted $n$-category on $C_{n}$. Thus, we must check that it satisfies the conditions stated in Definition~\ref{DnCatSS}.

We begin by proving that, for $k\in n$, $(C_{n},(\#^{k},\mathrm{sc}^{k}, \mathrm{tg}^{k}))$ is a single-sorted category. To this end, we consider the different items stated in Definition~\ref{DCat}.

\textsf{(C1)} For every $f\in C_{n}$, the following chain of equalities holds
\allowdisplaybreaks
\begin{align*}
\mathrm{sc}^{k}\left(
\mathrm{sc}^{k}\left(
f
\right)\right)
&=
\mathrm{i}^{(n,k)}\left(
\mathrm{sc}^{(k,n)}\left(
\mathrm{i}^{(n,k)}\left(
\mathrm{sc}^{(k,n)}\left(
f
\right)\right)\right)\right)
\tag{1}
\\&=
\mathrm{i}^{(n,k)}\left(
\mathrm{sc}^{(k,n)}\left(
f
\right)\right)
\tag{2}
\\&=
\mathrm{sc}^{k}\left(
f
\right).
\tag{3}
\end{align*}

The first equality unravels the definition of the unary operation $\mathrm{sc}^{k}$ introduced in Definition~\ref{DUnif}; the second equality follows from the fact that,  by item~(MS5) in Definition~\ref{DnCatMS}, $\mathrm{sc}^{(k,n)}\circ
\mathrm{i}^{(n,k)}=\mathrm{Id}_{C_{k}}$; finally, the last equality recovers the definition of the unary operation $\mathrm{sc}^{k}$ introduced in Definition~\ref{DUnif}.

By a similar argument we also have that
\begin{align*}
\mathrm{sc}^{k}\circ\mathrm{tg}^{k}
&=\mathrm{tg}^{k},
&
\mathrm{tg}^{k}\circ\mathrm{sc}^{k}
&=\mathrm{sc}^{k},
&
\mathrm{tg}^{k}\circ\mathrm{tg}^{k}
&=\mathrm{tg}^{k}.
\end{align*}

\textsf{(C2)} For every $f,g\in C_{n}$ the following chain of equivalences holds
\allowdisplaybreaks
\begin{align*}
g\#^{k}f\mbox{ is defined }
&\Leftrightarrow
g\circ^{(k)}f\mbox{ is defined}
\tag{1}
\\&\Leftrightarrow
\mathrm{sc}^{(k,n)}(g)=\mathrm{tg}^{(k,n)}(f)
\tag{2}
\\&\Leftrightarrow
\mathrm{sc}^{k}(g)=\mathrm{tg}^{k}(f).
\tag{3}
\end{align*}

The first equivalence unravels the definition of the partial binary operation $\#^{k}$ introduced in Definition~\ref{DUnif}; the second equivalence follows from the fact that $\mathsf{C}$ is a many-sorted $n$-category, see Definition~\ref{DnCatMS}; finally, the last equivalence follows from Corollary~\ref{CUnifScTg}.

\textsf{(C3)} For every $f,g\in C_{n}$, if $\mathrm{sc}^{k}(g)=\mathrm{tg}^{k}(f)$, the following chain of equalities holds
\allowdisplaybreaks
\begin{align*}
\mathrm{sc}^{k}\left(
g\#^{k}f
\right)
&=
\mathrm{i}^{(n,k)}\left(
\mathrm{sc}^{(k,n)}\left(
g\circ^{(k)}f
\right)\right)
\tag{1}
\\&=
\mathrm{i}^{(n,k)}\left(
\mathrm{sc}^{(k,n)}\left(
f
\right)\right)
\tag{2}
\\&=
\mathrm{sc}^{k}\left(
f
\right).
\tag{3}
\end{align*}

The first equality unravels the definition of the unary operation $\mathrm{sc}^{k}$ and the partial binary operation $\#^{k}$ introduced in Definition~\ref{DUnif}; the second equality follows from the fact that, by item~(MS2) in Definition~\ref{DnCatMS}, $\mathrm{sc}^{(k,n)}(g\circ^{(k)}f)=\mathrm{sc}^{(k,n)}(f)$; finally, the last equality recovers the definition of the unary operation $\mathrm{sc}^{k}$ introduced in Definition~\ref{DUnif}.

By a similar argument we also have that
$$
\mathrm{tg}^{k}\left(
g\#^{k}f
\right)=\mathrm{tg}^{k}(g).
$$

\textsf{(C4)} For every $f\in C_{n}$, the following chain of equalities holds
\allowdisplaybreaks
\begin{align*}
f\#^{k}\mathrm{sc}^{k}(f)&=
f\circ^{(k)}\left(
\mathrm{i}^{(n,k)}\left(
\mathrm{sc}^{(k,n)}\left(
f
\right)\right)\right)
\tag{1}
\\&=
f.
\tag{2}
\end{align*}
The first equality unravels the definition of the unary operation $\mathrm{sc}^{k}$ and the partial binary operation $\#^{k}$ introduced in Definition~\ref{DUnif}; the second equality follows from item~(MS6) in Definition~\ref{DnCatMS}.

By a similar argument we also have that
$$
\mathrm{tg}^{k}(f)\#^{k}f=f.
$$

\textsf{(C5)} For every $f,g,h\in C_{n}$, if $\mathrm{sc}^{k}(h)=\mathrm{tg}^{k}(g)$ and $\mathrm{sc}^{k}(g)=\mathrm{tg}^{k}(f)$, then the following chain of equalities holds
\allowdisplaybreaks
\begin{align*}
h\#^{k}\left(
g\#^{k}f
\right)&=
h\circ^{(k)}\left(
g\circ^{(k)}f\right)
\tag{1}
\\&=
\left(
h\circ^{(k)}g
\right)\circ^{(k)}f
\tag{2}
\\&=
\left(
h\#^{k}g\right)\#^{k}f.
\tag{3}
\end{align*}

The first equality unravels the definition of the partial binary operation $\#^{k}$ introduced in Definition~\ref{DUnif}; the second equality follows from item~(MS7) in Definition~\ref{DnCatMS}; finally, the last equality  recovers the definition of the partial binary operation $\#^{k}$ introduced in Definition~\ref{DUnif}.

It follows that, for every $k\in n$, $(C_{n}, \#^{k},\mathrm{sc}^{k}, \mathrm{tg}^{k})$ is a single-sorted category.

According to Definition~\ref{DnCatSS}, now we have to prove that, for every $j,k\in n$ with $j<k$, the ordered pair $(C_{n},(\#^{i},\mathrm{sc}^{i},\mathrm{tg}^{i})_{i\in\{j,k\}})$ is a $2$-category.
To this end, we consider the different items stated in Definition~\ref{DnCatSS}.

\textsf{(2C1)} For every $f\in C_{n}$, the following chain of equalities holds
\allowdisplaybreaks
\begin{align*}
\mathrm{sc}^{j}\left(
\mathrm{sc}^{k}\left(
f
\right)\right)&=
\mathrm{i}^{(n,j)}\left(
\mathrm{sc}^{(j,n)}\left(
\mathrm{i}^{(n,k)}\left(
\mathrm{sc}^{(k,n)}\left(
f
\right)\right)\right)\right)
\tag{1}
\\&=
\mathrm{i}^{(n,j)}\left(
\mathrm{sc}^{(j,k)}\left(
\mathrm{sc}^{(k,n)}\left(
\mathrm{i}^{(n,k)}\left(
\mathrm{sc}^{(k,n)}\left(
f
\right)\right)\right)\right)\right)
\tag{2}
\\&=
\mathrm{i}^{(n,j)}\left(
\mathrm{sc}^{(j,k)}\left(
\mathrm{sc}^{(k,n)}\left(
f
\right)\right)\right)
\tag{3}
\\&=
\mathrm{i}^{(n,j)}\left(
\mathrm{sc}^{(j,n)}\left(
f
\right)\right)
\tag{4}
\\&=
\mathrm{sc}^{j}\left(
f\right).
\tag{5}
\end{align*}

The first equality unravels the definition of the unary operations $\mathrm{sc}^{j}$ and $\mathrm{sc}^{k}$ introduced in Definition~\ref{DUnif}; the second equality follows from the fact that, by item~(MS1) from Definition~\ref{DnCatMS}, $\mathrm{sc}^{(j,n)}= \mathrm{sc}^{(j,k)}\circ \mathrm{sc}^{(k,n)}$; the third equality follows from the fact that, by item~(MS5) from Definition~\ref{DnCatMS},  $\mathrm{sc}^{(k,n)}
\circ
\mathrm{i}^{(n,k)}=\mathrm{Id}_{C_{n}}$; the fourth equality follows from the fact that, by item~(MS1) from Definition~\ref{DnCatMS}, $\mathrm{sc}^{(j,k)}\circ\mathrm{sc}^{(k,n)}=\mathrm{sc}^{(j,n)}$; finally, the last equality recovers the definition of the unary operations $\mathrm{sc}^{j}$ introduced in Definition~\ref{DUnif}.

By a similar argument we also have that 
$$
\mathrm{sc}^{j}\circ \mathrm{tg}^{k}=\mathrm{sc}^{j}.
$$

On the other hand, for every $f\in C_{n}$, the following chain of equalities also holds
\allowdisplaybreaks
\begin{align*}
\mathrm{sc}^{k}\left(
\mathrm{sc}^{j}\left(
f
\right)\right)&=
\mathrm{i}^{(n,k)}\left(
\mathrm{sc}^{(k,n)}\left(
\mathrm{i}^{(n,j)}\left(
\mathrm{sc}^{(j,n)}\left(
f
\right)\right)\right)\right)
\tag{1}
\\&=
\mathrm{i}^{(n,k)}\left(
\mathrm{sc}^{(k,n)}\left(
\mathrm{i}^{(n,k)}\left(
\mathrm{i}^{(k,j)}\left(
\mathrm{sc}^{(j,n)}\left(
f
\right)\right)\right)\right)\right)
\tag{2}
\\&=
\mathrm{i}^{(n,k)}\left(
\mathrm{i}^{(k,j)}\left(
\mathrm{sc}^{(j,n)}\left(
f
\right)\right)\right)
\tag{3}
\\&=
\mathrm{i}^{(n,j)}\left(
\mathrm{sc}^{(j,n)}\left(
f
\right)\right)
\tag{4}
\\&=
\mathrm{sc}^{j}(f).
\tag{5}
\end{align*}

The first equality unravels the definition of the unary operations $\mathrm{sc}^{j}$ and $\mathrm{sc}^{k}$ introduced in Definition~\ref{DUnif}; the second equality follows from the fact that, by item~(MS4) from Definition~\ref{DnCatMS}, $\mathrm{i}^{(n,j)}=\mathrm{i}^{(n,k)}\circ\mathrm{i}^{(k,j)}$; the third equality follows from the fact that, by item~(MS5) from Definition~\ref{DnCatMS},  $\mathrm{sc}^{(k,n)}
\circ
\mathrm{i}^{(n,k)}=\mathrm{Id}_{C_{n}}$; the fourth equality follows from the fact that, by item~(MS4) from Definition~\ref{DnCatMS}, $\mathrm{i}^{(n,j)}=\mathrm{i}^{(n,k)}\circ\mathrm{i}^{(k,j)}$; finally, the last equality recovers the definition of the unary operations $\mathrm{sc}^{j}$ introduced in Definition~\ref{DUnif}.

By a similar argument we also have that
$$
\mathrm{tg}^{k}\circ \mathrm{tg}^{j}
=
\mathrm{tg}^{j}
=
\mathrm{tg}^{j}\circ\mathrm{tg}^{k}
=
\mathrm{tg}^{j}\circ\mathrm{sc}^{k}.
$$

\textsf{(2C2)} For every $f,g\in C_{n}$, if $\mathrm{sc}^{j}(g)=\mathrm{tg}^{j}(f)$, the following chain of equalities holds
\allowdisplaybreaks
\begin{align*}
\mathrm{sc}^{k}\left(
g\#^{j}f
\right)&=
\mathrm{i}^{(n,k)}\left(
\mathrm{sc}^{(k,n)}\left(
g\circ^{(j)}f
\right)\right)
\tag{1}
\\&=
\mathrm{i}^{(n,k)}\left(
\mathrm{sc}^{(k,n)}\left(
g
\right)
\circ^{(j)}
\mathrm{sc}^{(k,n)}\left(
f
\right)\right)
\tag{2}
\\&=
\left(\left(\mathrm{i}^{(n,k)}\left(
\mathrm{sc}^{(k,n)}\right)\left(g
\right)\right)\right)
\circ^{(j)}
\left(\left(\mathrm{i}^{(n,k)}\left(
\mathrm{sc}^{(k,n)}\right)\left(
f
\right)\right)\right)
\tag{3}
\\&=
\mathrm{sc}^{k}(g)
\#^{j}
\mathrm{sc}^{k}(f).
\tag{4}
\end{align*}

The first equality unravels the definition of the unary operation $\mathrm{sc}^{k}$ and  the partial binary operation $\#^{j}$ introduced in Definition~\ref{DUnif}; the second equality follows from the fact that, by item~(MS3) from Definition~\ref{DnCatMS}, $\mathrm{sc}^{(k,n)}(
g\circ^{(j)}f
)= 
\mathrm{sc}^{(k,n)}(g)
\circ^{(j)}
\mathrm{sc}^{(k,n)}(f)$; the third equality follows from the fact that, by item~(MS8) from Definition~\ref{DnCatMS}, 
\begin{multline*}\mathrm{i}^{(n,k)}\left(
\mathrm{sc}^{(k,n)}\left(
g
\right)
\circ^{(j)}
\mathrm{sc}^{(k,n)}\left(
f
\right)\right)
\\=
\left(\left(
\mathrm{i}^{(n,k)}\left(
\mathrm{sc}^{(k,n)}\right)\left(
g
\right)\right)\right)
\circ^{(j)}
\left(\left(
\mathrm{i}^{(n,k)}\left(
\mathrm{sc}^{(k,n)}\right)\left(
f
\right)\right)\right);
\end{multline*}
finally, the last equality recovers the definition of the unary operation $\mathrm{sc}^{k}$ and  the partial binary operation $\#^{j}$ introduced in Definition~\ref{DUnif}.

By a similar argument we also have that
$$
\mathrm{tg}^{k}\left(
g\#^{j}f\right)
=
\mathrm{tg}^{k}(g)
\#^{j}
\mathrm{tg}^{k}(f).
$$

\textsf{(2C3)} For every $f,f',g,g'\in C_{n}$ if
\begin{align*}
\mathrm{sc}^{k}(g')&=\mathrm{tg}^{k}(g);
&
\mathrm{sc}^{k}(f')&=\mathrm{tg}^{k}(f);
\\
\mathrm{sc}^{j}(g')&=\mathrm{tg}^{j}(f');
&
\mathrm{sc}^{j}(g)&=\mathrm{tg}^{j}(f),
\end{align*}
then the following chain of equalities holds
\allowdisplaybreaks
\begin{align*}
\left(
g'\#^{j}f'
\right)\#^{k}\left(
g\#^{j} f\right)&=
\left(
g'\circ^{(j)}f'
\right)\circ^{(k)}
\left(
g\#^{(j)} f
\right)
\tag{1}
\\&=
\left(
g'\circ^{(k)}g\right)
\circ^{(j)}
\left(
f'\circ^{(k)}f
\right)
\tag{2}
\\&=
\left(
g'\#^{k}g
\right)
\circ^{j}\left(
f'\#^{k}f
\right).
\tag{3}
\end{align*}

The first equality unravels the definition of  the partial binary operations $\#^{j}$ and $\#^{k}$ introduced in Definition~\ref{DUnif}; the second equality follows from item~(MS9) in Definition~\ref{DnCatMS}; finally, the last equality recovers the definition of  the partial binary operations $\#^{j}$ and $\#^{k}$ introduced in Definition~\ref{DUnif}.

All in all, we can affirm that $\mathrm{Uf}^{(n)}(\mathsf{C})$ is a single-sorted $n$-category.
\end{proof}

In the following definition we introduce the notion of unification of a many-sorted $n$-functor.

\begin{definition}\label{DUnifFun} 
Let $\mathsf{C} = ((C_{k})_{k\in n+1},\zeta)$ with 
$
\zeta = (\zeta_{j,k})_{(j,k)\in \coprod_{k\in n+1}k}
$ 
and, for every $(j,k)\in \coprod_{k\in n+1}k$,
$$
\zeta_{j,k} =  \left(
\circ^{(j)},\mathrm{sc}^{(j,k)},\mathrm{tg}^{(j,k)},\mathrm{i}^{(k,j)}
\right), \text{ and }
$$
$\mathsf{C}' = ((C'_{k})_{k\in n+1},\zeta')$ with 
$
\zeta' = (\zeta'_{j,k})_{(j,k)\in \coprod_{k\in n+1}k}
$
and, for every $(j,k)\in \coprod_{k\in n+1}k$,
$$
\zeta'_{j,k} =  \left(
\circ'^{(j)},\mathrm{sc}'^{(j,k)},\mathrm{tg}'^{(j,k)},\mathrm{i}'^{(k,j)}
\right) 
$$
two many-sorted $n$-categories, $F=(F_{k})_{k\in n+1}\colon \mathsf{C}\mor\mathsf{C}'$ a many-sorted $n$-functor from $\mathsf{C}$ to $\mathsf{C}'$, and
\allowdisplaybreaks
\begin{align*}
\mathrm{Uf}^{(n)}\left(
\mathsf{C}
\right)&=\left(
C_{n}, \left(
\#^{k},\mathrm{sc}^{k},\mathrm{tg}^{k}
\right)_{k\in n}\right); 
&
\mathrm{Uf}^{(n)}\left(
\mathsf{C}'
\right)&=
\left(
C'_{n}, \left(
\#'^{k},\mathrm{sc}'^{k},\mathrm{tg}'^{k}
\right)_{k\in n}\right)
\end{align*}
the unifications corresponding to $\mathsf{C}$ and $\mathsf{C}'$, respectively. 

We will call the mapping $F_{n}$ the \emph{unification of $F$} and we will denote it by 
$\mathrm{Uf}^{(n)}(F)$. Thus $F_{n}$ is a mapping from $\mathrm{Uf}^{(n)}(\mathsf{C})$ to $\mathrm{Uf}^{(n)}(\mathsf{C}')$.
\end{definition}

\begin{proposition}\label{PUnifFun} 
Let $\mathsf{C} = ((C_{k})_{k\in n+1},\zeta)$ with 
$
\zeta = (\zeta_{j,k})_{(j,k)\in \coprod_{k\in n+1}k}
$ 
and, for every $(j,k)\in \coprod_{k\in n+1}k$,
$$
\zeta_{j,k} =  \left(
\circ^{(j)},\mathrm{sc}^{(j,k)},\mathrm{tg}^{(j,k)},\mathrm{i}^{(k,j)}
\right), \text{ and }
$$
$\mathsf{C}' = ((C'_{k})_{k\in n+1},\zeta')$ with 
$
\zeta' = (\zeta'_{j,k})_{(j,k)\in \coprod_{k\in n+1}k}
$ 
and, for every $(j,k)\in \coprod_{k\in n+1}k$,
$$
\zeta'_{j,k} = \left(\circ'^{(j)},\mathrm{sc}'^{(j,k)},\mathrm{tg}'^{(j,k)},\mathrm{i}'^{(k,j)}
\right) 
$$
two many-sorted $n$-categories and $F\colon \mathsf{C}\mor\mathsf{C}'$ a many-sorted $n$-functor from $\mathsf{C}$ to $\mathsf{C}'$. Then $\mathrm{Uf}^{(n)}(F)$  is a single-sorted $n$-functor from $\mathrm{Uf}^{(n)}(\mathsf{C})$ to $\mathrm{Uf}^{(n)}(\mathsf{C}')$.
\end{proposition}

\begin{proof}
It suffices to prove that the single-sorted mapping $\mathrm{Uf}^{(n)}(F)$ described in Definition~\ref{DUnifFun} satisfies the conditions of Definition~\ref{DnCatSS}, thus defining a single-sorted $n$-functor from $\mathrm{Uf}^{(n)}(\mathsf{C})$ to $\mathrm{Uf}^{(n)}(\mathsf{C}')$.

Let $k$ be a natural number in $n$.

\textsf{(Ci)} If $f\in C_{n}$, then the following chain of equalities holds
\allowdisplaybreaks
\begin{align*}
\mathrm{Uf}^{(n)}(F)\left(
\mathrm{sc}^{k}\left(
f
\right)\right)&=
F_{n}\left(
\mathrm{i}^{(n,k)}\left(
\mathrm{sc}^{(k,n)}\left(
f
\right)\right)\right)
\tag{1}
\\&=
\mathrm{i}'^{(n,k)}\left(
F_{k}\left(
\mathrm{sc}^{(k,n)}\left(
f
\right)\right)\right)
\tag{2}
\\&=
\mathrm{i}'^{(n,k)}\left(
\mathrm{sc}'^{(k,n)}\left(
F_{n}(f)
\right)\right)
\tag{3}
\\&=
\mathrm{sc}'^{k}\left(
\mathrm{Uf}^{(n)}(F)\left(
f
\right)\right).
\tag{4}
\end{align*}

The first equality unravels the definition of the unary operation $\mathrm{sc}^{k}$ introduced in Definition~\ref{DUnif} and the definition of the unification of $F$ introduced in Definition~\ref{DUnifFun}; the second equality follows from item~(MSii) in Definition~\ref{DnCatMS}; the third equality follows from item~(MSi) in Definition~\ref{DnCatMS}; finally, the last equality recovers the definition of the unary operation $\mathrm{sc}'^{k}$ introduced in Definition~\ref{DUnif} and the definition of the unification of $F$ introduced in Definition~\ref{DUnifFun}.

In an exactly similar way we also have that 
$$
\mathrm{Uf}^{(n)}(F)\left(
\mathrm{tg}^{k}\left(
f
\right)\right)=
\mathrm{tg}'^{k}\left(
\mathrm{Uf}^{(n)}(F)\left(
f
\right)\right).
$$

\textsf{(Cii)} If $f,g\in C_{n}$ are such that $\mathrm{sc}^{k}(g)=\mathrm{tg}^{k}(f)$, then the following chain of equalities holds
\allowdisplaybreaks
\begin{align*}
\mathrm{Uf}^{(n)}(F)\left(
g\#^{k}f
\right)&=
F_{n}\left(
g\circ^{(k)}f
\right)
\tag{1}
\\&=
F_{n}(g)\circ'^{(k)}F_{n}(f)
\tag{2}
\\&=
\mathrm{Uf}^{(n)}(F)(g)\#'^{k}\mathrm{Uf}^{(n)}(F)(f).
\tag{3}
\end{align*}

The first equality unravels the definition of the partial binary operation $\#^{k}$ introduced in Definition~\ref{DUnif} and the definition of the unification of $F$ introduced in Definition~\ref{DUnifFun};  the second equality follows from item~(MSiii) in Definition~\ref{DnCatMS}; finally, the last equality recovers the definition of the partial binary operation $\#^{k}$ introduced in Definition~\ref{DUnif} and the definition of the unification of $F$ introduced in Definition~\ref{DUnifFun}.

Therefore, $\mathrm{Uf}^{(n)}(F)$ is a single-sorted $n$-functor from $\mathrm{Uf}^{(n)}(\mathsf{C})$ to $\mathrm{Uf}^{(n)}(\mathsf{C}')$.
\end{proof}

We are now in position to prove that the unification construction is functorial.

\begin{proposition}\label{PUnifFunctor} 
$\mathrm{Uf}^{(n)}\colon \mathsf{nCat}^{\mathrm{MS}}\mor \mathsf{nCat}$ is a functor.
\end{proposition}
\begin{proof}
Taking into account Propositions~\ref{PUnif} and~\ref{PUnifFun}, it suffices to prove that $\mathrm{Uf}^{(n)}$ preserves many-sorted identity $n$-functors and compositions of many-sorted $n$-functors.

Let $\mathsf{C}$ be a many-sorted $n$-category and consider the many-sorted $n$-functor $\mathrm{Id}_{\mathsf{C}}\colon\mathsf{C}\mor\mathsf{C}$, where $\mathrm{Id}_{\mathsf{C}}=(\mathrm{Id}_{C_{k}})_{k\in n+1}$. The following chain of equalities holds
\begin{align*}
\mathrm{Uf}^{(n)}\left(
\mathrm{Id}_{\mathsf{C}}
\right)&=
\mathrm{Id}_{C_{n}}
\tag{1}
\\&=\mathrm{Id}_{\mathrm{Uf}^{(n)}(
\mathsf{C})}.
\tag{2}
\end{align*}

The first equality unravels Definition~\ref{DUnifFun} on the unification of $\mathrm{Id}_{\mathsf{C}}$; the second equality unravels Definition~\ref{DUnif} on the unification of $\mathsf{C}$. Note that, by Definition~\ref{DnCatSS}, this mapping corresponds to the single-sorted identity $n$-functor on $\mathrm{Uf}^{(n)}(\mathsf{C})$.

Thus, $\mathrm{Uf}^{(n)}$ preserves identities.

Now let $F\colon\mathsf{C}\mor \mathsf{C}'$ and $G\colon\mathsf{C}'\mor\mathsf{C}''$ be two many-sorted $n$-functors between many-sorted $n$-categories, where $F=(F_{k})_{k\in n+1}$, and $G=(G_{k})_{k\in n+1}$. Then we want to check that 
$$
\mathrm{Uf}^{(n)}\left(
G\circ F
\right)=\mathrm{Uf}^{(n)}(G)\circ\mathrm{Uf}^{(n)}(F).
$$

The following chain of equalities holds
\allowdisplaybreaks
\begin{align*}
\mathrm{Uf}^{(n)}\left(
G\circ F
\right)&=\mathrm{Uf}^{(n)}
\left(\left(
G_{k}\circ F_{k}
\right)_{k\in n+1}
\right)
\tag{1}
\\&=G_{n}\circ F_{n}
\tag{2}
\\&=\mathrm{Uf}^{(n)}(G)\circ\mathrm{Uf}^{(n)}(F).
\tag{3}
\end{align*}

The first equality unravels the definition of composition of many-sorted $n$-functors from Definition~\ref{DnCatMS}; the second equality unravels Definition~\ref{DUnifFun} on the unification of $G\circ F$; the third equality recovers Definition~\ref{DUnifFun} on the unification of, respectively, $G$ and $F$. Note that, by Definition~\ref{DnCatSS}, this mapping corresponds to the composition of the single-sorted identity $n$-functor $\mathrm{Uf}^{(n)}(G)$ and $\mathrm{Uf}^{(n)}(F)$.

Thus, $\mathrm{Uf}^{(n)}$ preserves composition.

This completes the proof.
\end{proof}

\section{
\texorpdfstring
{On the relationship between \(n\)-graduation and \(n\)-unification}
{On the relationship between n-graduation and n-unification}
}

In this section, we investigate the relation between the functors $\mathrm{Gd}^{(n)}$ and $\mathrm{Uf}^{(n)}$.
\begin{center}
\begin{tikzpicture}
[ACliment/.style={-{To [angle'=45, length=5.75pt, width=4pt, round]}}
, scale=1, 
AClimentD/.style={dogble eqgal sign distance,
-implies
}
]

\node[] (0R) at (0.5,0) [color=white] {};
\node[] (1L) at (2.3,0)  [color=white] {};

\node[] (C0) at (0,0) [] {$\mathsf{nCat}$};
\node[] (C1) at (3,0) [] {$\mathsf{nCat}^{\mathrm{MS}}$};

\draw[ACliment] (0R.north east) to node [above] {$\mathrm{Gd}^{(n)}$} (1L.north west);
\draw[ACliment] (1L.south west) to node [below] {$\mathrm{Uf}^{(n)}$} (0R.south east);
\end{tikzpicture}
\end{center}

We start by proving that the process of graduation followed by that of unification, does not provide anything new.

\begin{proposition}\label{PUnifGrad} 
$\mathrm{Uf}^{(n)}\circ\mathrm{Gd}^{(n)}=\mathrm{Id}_{\mathsf{nCat}}$.
\end{proposition}

\begin{proof}
Let $\mathsf{C}=(C,(\xi_{k})_{k\in n})$ be a single-sorted $n$-category, with $(\xi_{k})_{k\in n} = (\#^{k},\mathrm{sc}^{k},\mathrm{tg}^{k})_{k\in n}$, and 
$$
\mathrm{Gd}^{(n)}\left(
\mathsf{C}
\right) = 
\left(\left(
C_{k}
\right)_{k\in n+1}, 
\left(
\circ^{(j)},\mathrm{sc}^{(j,k)},\mathrm{tg}^{(j,k)},\mathrm{i}^{(k,j)}
\right)_{(j,k)\in \coprod_{k\in n+1}k}
\right)
$$
the graduation of $\mathsf{C}$. Then, according to Definition~\ref{DGrad}, for every $k\in n$,
$$
C_{k}=\left\lbrace f\in C\mid f=\mathrm{sc}^{k}(f) \right\rbrace
$$
and that $C_{n}$ is precisely $C$, the underlying set of $\mathsf{C}$.

Moreover, for every $0\leq j<k\leq n$, the structural operations of the graduation are given by
\allowdisplaybreaks
\begin{align*}
\circ^{(j)}&=\#^{j}\bigr|_{C_{k}\times C_{k}}^{C_{k}};
&
\mathrm{i}^{(k,j)}&=\mathrm{in}_{C_{j},C_{k}};
\\
\mathrm{sc}^{(j,k)}&=\mathrm{sc}^{j}\bigr|_{C_{k}}^{C_{j}};
&
\mathrm{tg}^{(j,k)}&=\mathrm{tg}^{j}\bigr|_{C_{k}}^{C_{j}}.
\end{align*}

Now, let us consider the single-sorted $n$-category given by the unification of the graduation of $\mathsf{C}$, i.e.,
$$
\mathrm{Uf}^{(n)}\left(
\mathrm{Gd}^{(n)}\left(
\mathsf{C}
\right)\right)=
\left(
C',
\left(
\#'^{k},\mathrm{sc}'^{k},\mathrm{tg}'^{k}
\right)_{k\in n}\right).
$$

Then, according to Definition~\ref{DUnif}, we have that $C'$ is $C_{n}$, which, as we already know, is equal to $C$, the underlying set of $\mathsf{C}$. Moreover, for every $k\in n$, the structural operations of the unification are given by
\allowdisplaybreaks
\begin{align*}
\#'^{k}&=\circ^{(k)};
&
\mathrm{sc}'^{(k)}&=\mathrm{i}^{(n,k)}\circ\mathrm{sc}^{(k,n)};
&
\mathrm{tg}'^{(k)}&=\mathrm{i}^{(n,k)}\circ\mathrm{tg}^{(k,n)}.
\end{align*}

If, for every $k\in n$, we replace equals by equals in the structural operations, then we have that
\allowdisplaybreaks
\begin{align*}
\#'^{k}&=\circ^{(k)}
=\#^{k}\bigr|_{C_{n}\times C_{n}}^{C_{n}}
=\#^{k}\bigr|_{C\times C}^{C}
=\#^{k};
\\
\mathrm{sc}'^{(k)}&=\mathrm{i}^{(n,k)}\circ\mathrm{sc}^{(k,n)}
=\mathrm{in}_{C_{k},C_{n}}\circ\mathrm{sc}^{k}\bigr|_{C_{n}}^{C_{n}}
=\mathrm{sc}^{k}\bigr|_{C}^{C}
=\mathrm{sc}^{k};
\\
\mathrm{tg}'^{(k)}&=\mathrm{i}^{(n,k)}\circ\mathrm{tg}^{(k,n)}
=\mathrm{in}_{C_{k},C_{n}}\circ\mathrm{tg}^{k}\bigr|_{C_{n}}^{C_{n}}
=\mathrm{tg}^{k}\bigr|_{C}^{C}
=\mathrm{tg}^{k}.
\end{align*}

That is, $\mathrm{Uf}^{(n)}(\mathrm{Gd}^{(n)}(\mathsf{C}))=\mathsf{C}$.

Let $F\colon\mathsf{C}\mor\mathsf{C}'$ be a single-sorted $n$-functor between single-sorted $n$-categories and $\mathrm{Gd}^{(n)}(F)\colon\mathrm{Gd}^{(n)}(\mathsf{C})\mor\mathrm{Gd}^{(n)}(\mathsf{C}')$ the graduation of $F$. Then, according to Definition~\ref{DGradFun}, we have that $\mathrm{Gd}^{(n)}(F) = (F_{k})_{k\in n+1} = (F\bigr|_{C_{k}}^{C'_{k}})_{k\in n+1}$. Let us recall that, by Definition~\ref{DGrad}, the set $C_{n}$ is equal to $C$, the underlying set of $\mathsf{C}$, and that the set $C'_{n}$ is equal to $C'$, the underlying set of $\mathsf{C}'$.  

Now, let us consider the unification of the graduation of $F$, i.e., $\mathrm{Uf}^{(n)}(\mathrm{Gd}^{(n)}(F))$, which, according to Definition~\ref{DUnifFun}, is given by 
$$
\mathrm{Uf}^{(n)}\left(
\mathrm{Gd}^{(n)}\left(
F
\right)\right)
=\mathrm{Uf}^{(n)}
\left(\left(F\bigr|_{C_{k}}^{C'_{k}}
\right)_{k\in n+1}\right)=F\bigr|_{C_{n}}^{C'_{n}}
=F\bigr|_{C}^{C'}=F.
$$
Thus, $\mathrm{Uf}^{(n)}(\mathrm{Gd}^{(n)}(F))=F$.

This completes the proof.
\end{proof}

In contrast to the previous result, the process of unification followed by that of graduation does not give rise to the identity functor on the category of many-sorted $n$-categories. However, we will show that both functors are naturally isomorphic.

\begin{proposition}\label{PGradUnif} 
$\mathrm{Gd}^{(n)}\circ\mathrm{Uf}^{(n)}\cong \mathrm{Id}_{\mathsf{nCat}^{\mathrm{MS}}}$.
\end{proposition}

\begin{proof} 
Let $\mathsf{C} = ((C_{k})_{k\in n+1},\zeta)$ be a many-sorted $n$-category with 
$ 
\zeta = (\zeta_{j,k})_{(j,k)\in \coprod_{k\in n+1}k}
$ 
and, for every $(j,k)\in \coprod_{k\in n+1}k$,
$$
\zeta_{j,k} =  \left(
\circ^{(j)},\mathrm{sc}^{(j,k)},\mathrm{tg}^{(j,k)},\mathrm{i}^{(k,j)}
\right),
$$
and
$
\mathrm{Uf}^{(n)}(\mathsf{C})
=
(C, (\#^{k},\mathrm{sc}^{k},\mathrm{tg}^{k})_{k\in n}) 
$
the unification of $\mathsf{C}$. Then, according to Definition~\ref{DUnif}, $C=C_{n}$ and, for every $k\in n$, the structural operations in the unification are given by 
\allowdisplaybreaks
\begin{align*}
\#^{k}&=\circ^{(k)},
&
\mathrm{sc}^{k}&=\mathrm{i}^{(n,k)}\circ\mathrm{sc}^{(k,n)},
&
\mathrm{tg}^{k}&=\mathrm{i}^{(n,k)}\circ\mathrm{tg}^{(k,n)}.
\end{align*}

Now, let us consider the many-sorted $n$-category given by the graduation of the unification of $\mathsf{C}$, i.e.,
$$
\mathrm{Gd}^{(n)}\left(
\mathrm{Uf}^{(n)}\left(
\mathsf{C}
\right)\right)
=
\left(\left(
C'_{k}
\right)_{k\in n+1}, 
\left(
\circ'^{(j)},\mathrm{sc}'^{(j,k)},\mathrm{tg}'^{(j,k)},\mathrm{i}'^{(k,j)}
\right)_{(j,k)\in \coprod_{k\in n+1}k}
\right).
$$

Then, according to Definition~\ref{DGrad}, we have that, for every $k\in n$,
$$
C'_{k}=
\left\lbrace f\in C\mid f=\mathrm{sc}^{k}(f)\right\rbrace,
$$
moreover, $C'_{n}=C$, and for every $0\leq j< k\leq n$, the structural operations of the graduation are given by
\allowdisplaybreaks
\begin{align*}
\circ'^{( j)}&=\#^{j}\bigr|_{C'_{k}\times C'_{k}}^{C'_{k}},
&
\mathrm{i}'^{(k,j)}&=\mathrm{in}_{C'_{j},C'_{k}},
\\
\mathrm{sc}'^{(j,k)}&=\mathrm{sc}^{j}\bigr|_{C'_{k}}^{C'_{j}},
&
\mathrm{tg}'^{(j,k)}&=\mathrm{tg}^{j}\bigr|_{C'_{k}}^{C'_{j}}.
\end{align*}
 
If we replace equals by equals in the different sets and structural operations, then,  for every $k\in n$, we have that
$$
C'_{k}=
\left\lbrace
f\in C\mid f=\mathrm{sc}^{k}(f)
\right\rbrace
=
\left\lbrace
f\in C_{n}\mid f=\mathrm{i}^{(n,k)}\left(
\mathrm{sc}^{(k,n)}\left(
f
\right)\right)\right\rbrace,
$$
moreover $C'_{n}=C=C_{n}$, and, for every $0\leq j< k\leq n$, we have that
\allowdisplaybreaks
\begin{align*}
\circ'^{( j)}&=\#^{j}\bigr|_{C'_{k}\times C'_{k}}^{C'_{k}}=\circ^{(j)}\bigr|_{C'_{k}\times C'_{k}}^{C'_{k}},
&
\mathrm{i}'^{(k,j)}&=\mathrm{in}^{C'_{j},C'_{k}},
\\
\mathrm{sc}'^{(j,k)}&=\mathrm{sc}^{j}\bigr|_{C'_{k}}^{C'_{j}}
=\left(
\mathrm{i}^{(n,j)}\circ\mathrm{sc}^{(j,n)}
\right)\bigr|_{C'_{k}}^{C'_{j}},
&
\mathrm{tg}'^{(j,k)}&=\mathrm{tg}^{j}\bigr|_{C'_{k}}^{C'_{j}}
=\left(
\mathrm{i}^{(n,j)}\circ\mathrm{tg}^{(j,n)}
\right)\bigr|_{C'_{k}}^{C'_{j}}.
\end{align*}
 
\begin{claim}\label{CGradUnifA}
For every $k\in n$, $\mathrm{i}^{(n,k)}[C_{k}]=C'_{k}$.
\end{claim}

Indeed, let $f$ be an element of $C_{k}$. Then the following chain of equalities holds
\allowdisplaybreaks
\begin{align*}
\mathrm{sc}^{k}\left(
\mathrm{i}^{(n,k)}\left(
f
\right)\right)&=
\mathrm{i}^{(n,k)}\left(
\mathrm{sc}^{(k,n)}\left(
\mathrm{i}^{(n,k)}\left(
f
\right)\right)\right)
\tag{1}
\\&=
\mathrm{i}^{(n,k)}(f).
\tag{2}
\end{align*}

The first equality unravels the definition of the unary operation symbol $\mathrm{sc}^{k}$; the second equality follows from item~(MS4) in Definition~\ref{DnCatMS}.

Conversely, let $f$ be an element of $C'_{k}$, then, taking into account the definition of the set $C'_{k}$, we have that  $f=\mathrm{i}^{(n,k)}(\mathrm{sc}^{(k,n)}(f))$. Note that $\mathrm{sc}^{(k,n)}(f)$ is an element in $C_{k}$, hence $f\in \mathrm{i}^{(n,k)}[C_{k}]$.

All in all, we can affirm that $\mathrm{i}^{(n,k)}[C_{k}]=C'_{k}$. 

Claim~\ref{CGradUnifA} follows.

\begin{claim}\label{CGradUnifB}
For every $k\in n$, $\mathrm{sc}^{(k,n)}[C'_{k}]=C_{k}$.
\end{claim}

The following chain of equalities holds
\allowdisplaybreaks
\begin{align*}
\mathrm{sc}^{(k,n)}\left[
C'_{k}
\right]&=
\mathrm{sc}^{(k,n)}\left[
\mathrm{i}^{(n,k)}\left[
C_{k}
\right]\right]
\tag{1}
\\&=C_{k}.
\tag{2}
\end{align*}

The first equality follows from the fact that, by Claim~\ref{CGradUnifA}, for every $k\in n$, $\mathrm{i}^{(n,k)}[C_{k}]=C'_{k}$; finally, the last equality follows from the fact that, by item~(MS5) in Definition~\ref{DnCatMS} we have that $\mathrm{sc}^{(k,n)}\circ\mathrm{i}^{(n,k)}=\mathrm{Id}_{C_{k}}$.

All in all, we can affirm that $\mathrm{i}^{(n,k)}[C_{k}]=C'_{k}$. 

Claim~\ref{CGradUnifB} follows.

We let $\alpha^{\mathsf{C}}=(\alpha^{\mathsf{C}}_{k})_{k\in n+1}$ stand for the $n+1$-sorted mapping defined,   for every $k\in n+1$, as
$$
\alpha^{\mathsf{C}}_{k}=
\begin{cases}
\mathrm{i}^{(n,k)},&\mbox{if }k\in n;
\\
\mathrm{Id}_{C_{n}},&\mbox{if }k=n,
\end{cases}
$$
and represent it by the diagram of Figure~\ref{FAlpha}.

\begin{figure}
\begin{center}
\resizebox{\columnwidth}{!}{%
\begin{tikzpicture}
[ACliment/.style={-{To [angle'=45, length=5.75pt, width=4pt, round]}}
, scale=1, 
AClimentD/.style={dogble eqgal sign distance,
-implies
}
]

\node[] (0R) at (.2,0) [color=white] {$A^{A}_{A}$};
\node[] (1L) at (2.8,0)  [color=white] {$A^{A}_{A}$};
\node[] (1R) at (3.2,0) [color=white] {$A^{A}_{A}$};
\node[] (2L) at (5.8,0)  [color=white] {$A^{A}_{A}$};
\node[] (2R) at (6.2,0) [color=white] {$A^{A}_{A}$};
\node[] (3L) at (8.8,0)  [color=white] {$A^{A}_{A}$};
\node[] (3R) at (9.2,0) [color=white] {$A^{A}_{A}$};
\node[] (4L) at (11.8,0)  [color=white] {$A^{A}_{A}$};

\node[] (0Rp) at (.2,-2) [color=white] {$A^{A}_{A}$};
\node[] (1Lp) at (2.8,-2)  [color=white] {$A^{A}_{A}$};
\node[] (1Rp) at (3.2,-2) [color=white] {$A^{A}_{A}$};
\node[] (2Lp) at (5.8,-2)  [color=white] {$A^{A}_{A}$};
\node[] (2Rp) at (6.2,-2) [color=white] {$A^{A}_{A}$};
\node[] (3Lp) at (8.8,-2)  [color=white] {$A^{A}_{A}$};
\node[] (3Rp) at (9.2,-2) [color=white] {$A^{A}_{A}$};
\node[] (4Lp) at (11.8,-2)  [color=white] {$A^{A}_{A}$};

\node[] (C0) at (0,0) [] {$C_{0}$};
\node[] (C1) at (3,0) [] {$C_{1}$};
\node[] (C2) at (6,0) [] {$\cdots$};
\node[] (C3) at (9,0) [] {$C_{n-1}$};
\node[] (C4) at (12,0) [] {$C_{n}$};
\node[] (C0p) at (0,-2) [] {$C'_{0}$};
\node[] (C1p) at (3,-2) [] {$C'_{1}$};
\node[] (C2p) at (6,-2) [] {$\cdots$};
\node[] (C3p) at (9,-2) [] {$C'_{n-1}$};
\node[] (C4p) at (12,-2) [] {$C_{n}$};

\draw[ACliment] (0R) to node [midway, fill=white] {$\scriptstyle \mathrm{i}^{(1,0)}$} (1L);
\draw[ACliment] (1L.north west) to node [above] {$\scriptstyle \mathrm{sc}^{(0,1)}$} (0R.north east);
\draw[ACliment] (1L.south west) to node [below] {$\scriptstyle \mathrm{tg}^{(0,1)}$} (0R.south east);

\draw[ACliment] (1R) to node [midway, fill=white] {$\scriptstyle \mathrm{i}^{(2,1)}$} (2L);
\draw[ACliment] (2L.north west) to node [above] {$\scriptstyle \mathrm{sc}^{(1,2)}$} (1R.north east);
\draw[ACliment] (2L.south west) to node [below] {$\scriptstyle \mathrm{tg}^{(1,2)}$} (1R.south east);

\draw[ACliment] (2R) to node [midway, fill=white] {$\scriptstyle \mathrm{i}^{(n-1,n-2)}$} (3L);
\draw[ACliment] (3L.north west) to node [above] {$\scriptstyle \mathrm{sc}^{(n-2,n-1)}$} (2R.north east);
\draw[ACliment] (3L.south west) to node [below] {$\scriptstyle \mathrm{tg}^{(n-2,n-1)}$} (2R.south east);

\draw[ACliment] (3R) to node [midway, fill=white] {$\scriptstyle \mathrm{i}^{(n,n-1)}$} (4L);
\draw[ACliment] (4L.north west) to node [above] {$\scriptstyle \mathrm{sc}^{(n-1,n)}$} (3R.north east);
\draw[ACliment] (4L.south west) to node [below] {$\scriptstyle \mathrm{tg}^{(n-1,n)}$} (3R.south east);

\draw[ACliment] (0Rp) to node [midway, fill=white] {$\scriptstyle \mathrm{i}'^{(1,0)}$} (1Lp);
\draw[ACliment] (1Lp.north west) to node [above] {$\scriptstyle \mathrm{sc}'^{(0,1)}$} (0Rp.north east);
\draw[ACliment] (1Lp.south west) to node [below] {$\scriptstyle \mathrm{tg}'^{(0,1)}$} (0Rp.south east);

\draw[ACliment] (1Rp) to node [midway, fill=white] {$\scriptstyle \mathrm{i}'^{(2,1)}$} (2Lp);
\draw[ACliment] (2Lp.north west) to node [above] {$\scriptstyle \mathrm{sc}'^{(1,2)}$} (1Rp.north east);
\draw[ACliment] (2Lp.south west) to node [below] {$\scriptstyle \mathrm{tg}'^{(1,2)}$} (1Rp.south east);

\draw[ACliment] (2Rp) to node [midway, fill=white] {$\scriptstyle \mathrm{i}'^{(n-1,n-2)}$} (3Lp);
\draw[ACliment] (3Lp.north west) to node [above] {$\scriptstyle \mathrm{sc}'^{(n-2,n-1)}$} (2Rp.north east);
\draw[ACliment] (3Lp.south west) to node [below] {$\scriptstyle \mathrm{tg}'^{(n-2,n-1)}$} (2Rp.south east);

\draw[ACliment] (3Rp) to node [midway, fill=white] {$\scriptstyle \mathrm{i}'^{(n,n-1)}$} (4Lp);
\draw[ACliment] (4Lp.north west) to node [above] {$\scriptstyle \mathrm{sc}'^{(n-1,n)}$} (3Rp.north east);
\draw[ACliment] (4Lp.south west) to node [below] {$\scriptstyle \mathrm{tg}'^{(n-1,n)}$} (3Rp.south east);

\draw[ACliment, shorten <=0.2cm, shorten >=0.2cm] (C0) to node [left] {$\scriptstyle \mathrm{i}^{(n,0)}$} (C0p);
\draw[ACliment, shorten <=0.2cm, shorten >=0.2cm] (C1) to node [right] {$\scriptstyle \mathrm{i}^{(n,1)}$} (C1p);
\draw[ACliment, shorten <=0.2cm, shorten >=0.2cm] (C3) to node [right] {$\scriptstyle \mathrm{i}^{(n,n-1)}$} (C3p);
\draw[ACliment, shorten <=0.2cm, shorten >=0.2cm] (C4) to node [right] {$\scriptstyle \mathrm{Id}_{C_{n}}$} (C4p);
\end{tikzpicture}
 }
\end{center}
\caption{The $n+1$-sorted mapping $\alpha^{\mathsf{C}}$.}
\label{FAlpha}
\end{figure}

On the other hand, we let $\beta^{\mathsf{C}}=(\beta^{\mathsf{C}}_{k})_{k\in n+1}$ stand for the $n+1$-sorted mapping defined, for every $k\in n+1$, as   
$$
\beta^{\mathsf{C}}_{k}=
\begin{cases}
\mathrm{sc}^{(n,k)},&\mbox{if }k\in n;
\\
\mathrm{Id}_{C_{n}},&\mbox{if }k=n,
\end{cases}
$$
and represent it by the diagram in Figure~\ref{FBeta}.

\begin{figure}
\begin{center}
\resizebox{\columnwidth}{!}{%
\begin{tikzpicture}
[ACliment/.style={-{To [angle'=45, length=5.75pt, width=4pt, round]}}
, scale=1, 
AClimentD/.style={dogble eqgal sign distance,
-implies
}
]

\node[] (0R) at (.2,0) [color=white] {$A^{A}_{A}$};
\node[] (1L) at (2.8,0)  [color=white] {$A^{A}_{A}$};
\node[] (1R) at (3.2,0) [color=white] {$A^{A}_{A}$};
\node[] (2L) at (5.8,0)  [color=white] {$A^{A}_{A}$};
\node[] (2R) at (6.2,0) [color=white] {$A^{A}_{A}$};
\node[] (3L) at (8.8,0)  [color=white] {$A^{A}_{A}$};
\node[] (3R) at (9.2,0) [color=white] {$A^{A}_{A}$};
\node[] (4L) at (11.8,0)  [color=white] {$A^{A}_{A}$};

\node[] (0Rp) at (.2,-2) [color=white] {$A^{A}_{A}$};
\node[] (1Lp) at (2.8,-2)  [color=white] {$A^{A}_{A}$};
\node[] (1Rp) at (3.2,-2) [color=white] {$A^{A}_{A}$};
\node[] (2Lp) at (5.8,-2)  [color=white] {$A^{A}_{A}$};
\node[] (2Rp) at (6.2,-2) [color=white] {$A^{A}_{A}$};
\node[] (3Lp) at (8.8,-2)  [color=white] {$A^{A}_{A}$};
\node[] (3Rp) at (9.2,-2) [color=white] {$A^{A}_{A}$};
\node[] (4Lp) at (11.8,-2)  [color=white] {$A^{A}_{A}$};

\node[] (C0) at (0,0) [] {$C_{0}$};
\node[] (C1) at (3,0) [] {$C_{1}$};
\node[] (C2) at (6,0) [] {$\cdots$};
\node[] (C3) at (9,0) [] {$C_{n-1}$};
\node[] (C4) at (12,0) [] {$C_{n}$};
\node[] (C0p) at (0,-2) [] {$C'_{0}$};
\node[] (C1p) at (3,-2) [] {$C'_{1}$};
\node[] (C2p) at (6,-2) [] {$\cdots$};
\node[] (C3p) at (9,-2) [] {$C'_{n-1}$};
\node[] (C4p) at (12,-2) [] {$C_{n}$};

\draw[ACliment] (0R) to node [midway, fill=white] {$\scriptstyle \mathrm{i}^{(1,0)}$} (1L);
\draw[ACliment] (1L.north west) to node [above] {$\scriptstyle \mathrm{sc}^{(0,1)}$} (0R.north east);
\draw[ACliment] (1L.south west) to node [below] {$\scriptstyle \mathrm{tg}^{(0,1)}$} (0R.south east);

\draw[ACliment] (1R) to node [midway, fill=white] {$\scriptstyle \mathrm{i}^{(2,1)}$} (2L);
\draw[ACliment] (2L.north west) to node [above] {$\scriptstyle \mathrm{sc}^{(1,2)}$} (1R.north east);
\draw[ACliment] (2L.south west) to node [below] {$\scriptstyle \mathrm{tg}^{(1,2)}$} (1R.south east);

\draw[ACliment] (2R) to node [midway, fill=white] {$\scriptstyle \mathrm{i}^{(n-1,n-2)}$} (3L);
\draw[ACliment] (3L.north west) to node [above] {$\scriptstyle \mathrm{sc}^{(n-2,n-1)}$} (2R.north east);
\draw[ACliment] (3L.south west) to node [below] {$\scriptstyle \mathrm{tg}^{(n-2,n-1)}$} (2R.south east);

\draw[ACliment] (3R) to node [midway, fill=white] {$\scriptstyle \mathrm{i}^{(n,n-1)}$} (4L);
\draw[ACliment] (4L.north west) to node [above] {$\scriptstyle \mathrm{sc}^{(n-1,n)}$} (3R.north east);
\draw[ACliment] (4L.south west) to node [below] {$\scriptstyle \mathrm{tg}^{(n-1,n)}$} (3R.south east);

\draw[ACliment] (0Rp) to node [midway, fill=white] {$\scriptstyle \mathrm{i}'^{(1,0)}$} (1Lp);
\draw[ACliment] (1Lp.north west) to node [above] {$\scriptstyle \mathrm{sc}'^{(0,1)}$} (0Rp.north east);
\draw[ACliment] (1Lp.south west) to node [below] {$\scriptstyle \mathrm{tg}'^{(0,1)}$} (0Rp.south east);

\draw[ACliment] (1Rp) to node [midway, fill=white] {$\scriptstyle \mathrm{i}'^{(2,1)}$} (2Lp);
\draw[ACliment] (2Lp.north west) to node [above] {$\scriptstyle \mathrm{sc}'^{(1,2)}$} (1Rp.north east);
\draw[ACliment] (2Lp.south west) to node [below] {$\scriptstyle \mathrm{tg}'^{(1,2)}$} (1Rp.south east);

\draw[ACliment] (2Rp) to node [midway, fill=white] {$\scriptstyle \mathrm{i}'^{(n-1,n-2)}$} (3Lp);
\draw[ACliment] (3Lp.north west) to node [above] {$\scriptstyle \mathrm{sc}'^{(n-2,n-1)}$} (2Rp.north east);
\draw[ACliment] (3Lp.south west) to node [below] {$\scriptstyle \mathrm{tg}'^{(n-2,n-1)}$} (2Rp.south east);

\draw[ACliment] (3Rp) to node [midway, fill=white] {$\scriptstyle \mathrm{i}'^{(n,n-1)}$} (4Lp);
\draw[ACliment] (4Lp.north west) to node [above] {$\scriptstyle \mathrm{sc}'^{(n-1,n)}$} (3Rp.north east);
\draw[ACliment] (4Lp.south west) to node [below] {$\scriptstyle \mathrm{tg}'^{(n-1,n)}$} (3Rp.south east);

\draw[ACliment, shorten <=0.2cm, shorten >=0.2cm] (C0p) to node [left] {$\scriptstyle \mathrm{sc}^{(0,n)}$} (C0);
\draw[ACliment, shorten <=0.2cm, shorten >=0.2cm] (C1p) to node [right] {$\scriptstyle \mathrm{sc}^{(1,n)}$} (C1);
\draw[ACliment, shorten <=0.2cm, shorten >=0.2cm] (C3p) to node [right] {$\scriptstyle \mathrm{sc}^{(n-1,n)}$} (C3);
\draw[ACliment, shorten <=0.2cm, shorten >=0.2cm] (C4p) to node [right] {$\scriptstyle \mathrm{Id}_{C_{n}}$} (C4);
\end{tikzpicture}
}
\end{center}
\caption{The $n+1$-sorted mapping $\beta^{\mathsf{C}}$.}
\label{FBeta}
\end{figure}

Our next aim is to prove that, for every many-sorted $n$-category $\mathsf{C}$, $\alpha^{\mathsf{C}}$ and $\beta^{\mathsf{C}}$ are many-sorted $n$-functors.

\begin{claim}\label{CGradUnifC}
$\alpha^{\mathsf{C}}$ is a many-sorted $n$-functor from $\mathsf{C}$ to $\mathrm{Gd}^{(n)}(\mathrm{Uf}^{(n)}(\mathsf{C}))$.
\end{claim}

It suffices to prove that $\alpha^{\mathsf{C}}$ satisfies the conditions of Definition~\ref{DnCatMS}.

\textsf{(MSi)} For every $0\leq j< k \leq n$, if $f\in C_{k}$, the following chain of equalities holds
\allowdisplaybreaks
\begin{align*}
\alpha^{\mathsf{C}}_{j}\left(
\mathrm{sc}^{(j,k)}\left(
f
\right)\right)&=
\mathrm{i}^{(n,j)}\left(
\mathrm{sc}^{(j,k)}\left(
f
\right)\right)
\tag{1}
\\&=
\mathrm{i}^{(n,j)}\left(
\mathrm{sc}^{(j,k)}\left(
\mathrm{sc}^{(k,n)}\left(
\mathrm{i}^{(n,k)}\left(
f
\right)\right)\right)\right)
\tag{2}
\\&=
\mathrm{i}^{(n,j)}\left(
\mathrm{sc}^{(j,n)}\left(
\mathrm{i}^{(n,k)}\left(
f
\right)\right)\right)
\tag{3}
\\&=
\mathrm{sc}'^{(j,k)}\left(
\mathrm{i}^{(n,k)}\left(
f
\right)\right)
\tag{4}
\\&=
\mathrm{sc}'^{(j,k)}\left(
\alpha^{\mathsf{C}}_{k}\left(
f
\right)\right).
\tag{5}
\end{align*}

The first equality unravels the definition of the $j$-th component of $\alpha^{\mathsf{C}}$; the second equality follows from the fact that, by item~(MS5) in Definition~\ref{DnCatMS}, $\mathrm{sc}^{(k,n)}\circ\mathrm{i}^{(n,k)}=\mathrm{Id}_{C_{k}}$; the third equality follows from the fact that, by item~(MS1) in Definition~\ref{DnCatMS}, 
$\mathrm{sc}^{(j,k)}\circ\mathrm{sc}^{(k,n)}=\mathrm{sc}^{(j,n)}$; the fourth equality recovers the definition of the mapping $\mathrm{sc}'^{(j,k)}$ in $\mathrm{Gd}^{(n)}(\mathrm{Uf}^{(n)}(\mathsf{C}))$; finally, the last equality recovers the $k$-th component of $\alpha^{\mathsf{C}}$. 

In a similar fashion we have that 
$$
\alpha^{\mathsf{C}}_{j}\left(
\mathrm{tg}^{(j,k)}\left(
f
\right)\right)
=\mathrm{tg}'^{(j,k)}\left(
\alpha^{\mathsf{C}}_{k}\left(
f
\right)\right).
$$

\textsf{(MSii)} For every $0\leq  k <l< n$, if $f\in C_{k}$, the following chain of equalities holds
\allowdisplaybreaks
\begin{align*}
\alpha^{\mathsf{C}}_{l}\left(
\mathrm{i}^{(l,k)}\left(
f
\right)\right)&=
\mathrm{i}^{(n,l)}\left(
\mathrm{i}^{(l,k)}\left(
f
\right)\right)
\tag{1}
\\&=
\mathrm{i}^{(n,k)}\left(
f
\right)
\tag{2}
\\&=
\mathrm{in}^{C'_{k}}\left(
\mathrm{i}^{(n,k)}\left(
f
\right)\right)
\tag{3}
\\&=
\mathrm{i}'^{(l,k)}\left(
\mathrm{i}^{(n,k)}\left(
f
\right)\right)
\tag{4}
\\&=
\mathrm{i}'^{(l,k)}\left(
\alpha^{\mathsf{C}}_{k}\left(
f
\right)\right).
\tag{5}
\end{align*}

The first equality unravels the definition of the $l$-th component of $\alpha^{\mathsf{C}}$; the second equality follows from the fact that, by item~(MS4) in Definition~\ref{DnCatMS}, $\mathrm{i}^{(n,l)}\circ\mathrm{i}^{(l,k)}=\mathrm{i}^{(n,k)}$; the third equality applies the inclusion mapping $\mathrm{in}_{C'_{k}}\colon C'_{k}\mor C'_{l}$; the fourth equality recovers the definition of the mapping $\mathrm{i}'^{(l,k)}$ in $\mathrm{Gd}^{(n)}(\mathrm{Uf}^{(n)}(\mathsf{C}))$; finally, the last equality recovers the $k$-th component of $\alpha^{\mathsf{C}}$. 

For the particular case at the $n$-th sort, for every $0\leq  k < n$, the following chain of equalities holds
\allowdisplaybreaks
\begin{align*}
\alpha^{\mathsf{C}}_{n}\left(
\mathrm{i}^{(n,k)}\left(
f
\right)\right)&=
\mathrm{Id}_{C_{n}}\left(
\mathrm{i}^{(n,k)}\left(
f
\right)\right)
\tag{1}
\\&=
\mathrm{i}^{(n,k)}\left(
f
\right)
\tag{2}
\\&=
\mathrm{in}^{C'_{k}}\left(
\mathrm{i}^{(n,k)}\left(
f
\right)\right)
\tag{3}
\\&=
\mathrm{i}'^{(n,k)}\left(
\mathrm{i}^{(n,k)}\left(
f
\right)\right)
\tag{4}
\\&=
\mathrm{i}'^{(n,k)}\left(
\alpha^{\mathsf{C}}_{k}\left(
f
\right)\right).
\tag{5}
\end{align*}

The first equality unravels the definition of the $n$-th component of $\alpha^{\mathsf{C}}$; the second equality applies the identity mapping at $C_{n}$; the third equality applies the inclusion mapping $\mathrm{in}_{C'_{k}}\colon C'_{k}\mor C'_{n}$; the fourth equality recovers the definition of the mapping $\mathrm{i}'^{(n,k)}$ in $\mathrm{Gd}^{(n)}(\mathrm{Uf}^{(n)}(\mathsf{C}))$; finally, the last equality recovers the $k$-th component of $\alpha^{\mathsf{C}}$.

\textsf{(MSiii)} For every $0\leq j <k <n$, if $f,g\in C_{k}$ are such that $\mathrm{sc}^{(j,k)}(g)=\mathrm{tg}^{(j,k)}(f)$, then the following chain of equalities holds
\allowdisplaybreaks
\begin{align*}
\alpha^{\mathsf{C}}_{k}\left(
g\circ^{(j)}f
\right)&=
\mathrm{i}^{(n,k)}\left(
g\circ^{(j)}f
\right)
\tag{1}
\\&=
\mathrm{i}^{(n,k)}\left(
g
\right)
\circ^{(j)}
\mathrm{i}^{(n,k)}\left(
f
\right)
\tag{2}
\\&=
\mathrm{i}^{(n,k)}\left(
g
\right)\circ'^{(j)}\mathrm{i}^{(n,k)}\left(
f
\right)
\tag{3}
\\&=
\alpha^{\mathsf{C}}_{k}\left(
g
\right)\circ'^{(j)}
\alpha^{\mathsf{C}}_{k}\left(
f
\right).
\tag{4}
\end{align*} 

The first equality unravels the definition of the $k$-th component of $\alpha^{\mathsf{C}}$; the second equality follows from item~(MS8) in Definition~\ref{DnCatMS}; the third equality recovers the definition of the partial binary operation $\circ'^{(j)}$ in $\mathrm{Gd}^{(n)}(\mathrm{Uf}^{(n)}(\mathsf{C}))$; finally, the last equality recovers the $k$-th component of $\alpha^{\mathsf{C}}$.

For the particular case at the $n$-th sort, for every $0\leq  j < n$, the following chain of equalities holds
\allowdisplaybreaks
\begin{align*}
\alpha^{\mathsf{C}}_{n}\left(
g\circ^{(j)}f
\right)&=
\mathrm{Id}_{C_{n}}\left(
g\circ^{(j)}f
\right)
\tag{1}
\\&=
g\circ^{(j)}f
\tag{2}
\\&=
g\circ'^{(j)}f
\tag{3}
\\&=
\mathrm{Id}_{C_{n}}(g)\circ^{(j)}\mathrm{Id}_{C_{n}}(f)
\tag{4}
\\&=
\alpha^{\mathsf{C}}_{n}(g)\circ'^{(j)}\alpha^{\mathsf{C}}_{n}(f).
\tag{5}
\end{align*} 

The first equality unravels the definition of the $n$-th component of $\alpha^{\mathsf{C}}$; the second equality applies the identity mapping at $C_{n}$; the third equality recovers the definition of the partial binary operation $\circ'^{(j)}$ in $\mathrm{Gd}^{(n)}(\mathrm{Uf}^{(n)}(\mathsf{C}))$; the fourth equality recovers the identity mapping at $C_{n}$; finally, the last equality recovers the $n$-th component of $\alpha^{\mathsf{C}}$.

All in all, we can affirm that $\alpha^{\mathsf{C}}$ is a many-sorted $n$-functor from $\mathsf{C}$ to $\mathrm{Gd}^{(n)}(\mathrm{Uf}^{(n)}(\mathsf{C}))$,
$$\alpha^{\mathsf{C}}\colon \mathsf{C}\mor\mathrm{Gd}^{(n)}(\mathrm{Uf}^{(n)}(\mathsf{C})).$$

Claim~\ref{CGradUnifC} follows.

\begin{claim}\label{CGradUnifD}
$\beta^{\mathsf{C}}$ is a many-sorted $n$-functor from $\mathrm{Gd}^{(n)}(\mathrm{Uf}^{(n)}(\mathsf{C}))$ to $\mathsf{C}$.
\end{claim}

It suffices to prove that $\beta^{\mathsf{C}}$ satisfies the conditions of Definition~\ref{DnCatMS}.

\textsf{(MSi)} For every $0\leq j<k\leq n$, if $f\in C'_{k}$, the following chain of equalities holds
\allowdisplaybreaks
\begin{align*}
\beta_{j}^{\mathsf{C}}\left(
\mathrm{sc}'^{(j,k)}\left(
f
\right)\right)&=
\mathrm{sc}^{(j,n)}\left(
\mathrm{sc}'^{(j,k)}\left(
f
\right)\right)
\tag{1}
\\&=
\mathrm{sc}^{(j,n)}\left(
\mathrm{i}^{(n,j)}\left(
\mathrm{sc}^{(j,n)}\left(
f
\right)\right)\right)
\tag{2}
\\&=
\mathrm{sc}^{(j,n)}(f)
\tag{3}
\\&=
\mathrm{sc}^{(j,k)}\left(
\mathrm{sc}^{(k,n)}\left(
f
\right)\right)
\tag{4}
\\&=
\mathrm{sc}^{(j,k)}\left(
\beta^{\mathsf{C}}_{k}\left(
f
\right)\right).
\tag{5}
\end{align*}

The first equality unravels the definition of the $j$-th component of $\beta^{\mathsf{C}}$; the second equality unravels the definition of the mapping $\mathrm{sc}'^{(j,k)}$ in $\mathrm{Gd}^{(n)}(\mathrm{Uf}^{(n)}(\mathsf{C}))$; the third equality follows from the fact that, by item~(MS5) in Definition~\ref{DnCatMS}, $\mathrm{sc}^{(j,n)}\circ\mathrm{i}^{(n,j)}=\mathrm{Id}_{C_{j}}$; the fourth equality follows from the fact that, by item~(MS1) in Definition~\ref{DnCatMS}, $\mathrm{sc}^{(j,k)}\circ\mathrm{sc}^{(k,n)}=\mathrm{sc}^{(j,n)}$; finally, the last equality recovers the $k$-th component of $\beta^{\mathsf{C}}$.

By a similar argument we also have that 
$$
\beta^{\mathsf{C}}_{j}\left(
\mathrm{tg}'^{(j,k)}\left(
f
\right)\right)=\mathrm{tg}^{(j,k)}\left(
\beta^{\mathsf{C}}_{k}\left(
f
\right)\right).
$$

\textsf{(MSii)} For every $0\leq k<l<n$, if $f\in C'_{k}$, the following chain of equalities holds
\allowdisplaybreaks
\begin{align*}
\beta^{\mathsf{C}}_{l}\left(
\mathrm{i}'^{(l,k)}\left(
f
\right)\right)
&=
\mathrm{sc}^{(l,n)}\left(
\mathrm{i}'^{(l,k)}\left(
f
\right)\right)
\tag{1}
\\&=
\mathrm{sc}^{(l,n)}\left(
f
\right)
\tag{2}
\\&=
\mathrm{sc}^{(l,n)}\left(
\mathrm{i}^{(n,k)}\left(
\mathrm{sc}^{(k,n)}\left(
f
\right)\right)\right)
\tag{3}
\\&=
\mathrm{sc}^{(l,n)}\left(
\mathrm{i}^{(n,l)}\left(
\mathrm{i}^{(l,k)}\left(
\mathrm{sc}^{(k,n)}\left(
f
\right)\right)\right)\right)
\tag{4}
\\&=
\mathrm{i}^{(l,k)}\left(
\mathrm{sc}^{(k,n)}\left(
f
\right)\right)
\tag{5}
\\&=
\mathrm{i}^{(l,k)}\left(
\beta^{\mathsf{C}}_{k}\left(
f
\right)\right).
\tag{6}
\end{align*}

The first equality unravels the definition of the $l$-th component of $\beta^{\mathsf{C}}$; the second equality unravels the definition of the mapping $\mathrm{i}'^{(j,k)}$ in $\mathrm{Gd}^{(n)}(\mathrm{Uf}^{(n)}(\mathsf{C}))$. Let us recall that $\mathrm{i}'^{(l,k)}$ is just the inclusion mapping from $C'_{k}$ to $C'_{l}$;  the third equality follows from the fact that $f\in C'_{k}$, then taking into account the definition of the set $C'_{k}$, we have that  $f=\mathrm{i}^{(n,k)}(\mathrm{sc}^{(k,n)}(f))$; the fourth equality follows from the fact that, by item~(MS4) in Definition~\ref{DnCatMS}, $\mathrm{i}^{(n,l)}\circ\mathrm{i}^{(n,k)}=\mathrm{i}^{(n,k)}$; the fifth equality follows from the fact that, by item~(MS5) in Definition~\ref{DnCatMS}, $\mathrm{sc}^{(l,n)}\circ\mathrm{i}^{(n,l)}=\mathrm{Id}_{C_{l}}$; finally, the last equality recovers the $k$-th component of $\beta^{\mathsf{C}}$.

For the particular case at the $n$-th sort, for every $0\leq k<n$, the following chain of equalities holds.
\allowdisplaybreaks
\begin{align*}
\beta^{\mathsf{C}}_{n}\left(
\mathrm{i}'^{(n,k)}\left(
f
\right)\right)
&=
\mathrm{i}'^{(n,k)}\left(
f
\right)
\tag{1}
\\&=
f
\tag{2}
\\&=
\mathrm{i}^{(n,k)}\left(
\mathrm{sc}^{(k,n)}\left(
f
\right)\right)
\tag{3}
\\&=
\mathrm{i}^{(n,k)}\left(
\beta^{\mathsf{C}}_{k}\left(
f
\right)\right).
\tag{4}
\end{align*}

The first equality unravels the definition of the $n$-th component of $\beta^{\mathsf{C}}$. Let us recall that 
$\beta^{\mathsf{C}}_{n}=\mathrm{Id}_{C_{n}}$; the second equality unravels the definition of the mapping $\mathrm{i}'^{(n,k)}$ in $\mathrm{Gd}^{(n)}(\mathrm{Uf}^{(n)}(\mathsf{C}))$. Let us recall that $\mathrm{i}'^{(n,k)}$ is just the inclusion mapping from $C'_{k}$ to $C'_{n}$; the third equality follows from the fact that $f\in C'_{k}$, then taking into account the definition of the set $C'_{k}$, we have that  $f=\mathrm{i}^{(n,k)}(\mathrm{sc}^{(k,n)}(f))$; finally, the last equality recovers the $k$-th component of $\beta^{\mathsf{C}}$.

\textsf{(MSiii)} For every $0\leq j<k<n$, if $f,g\in C'_{k}$ are such that $\mathrm{sc}'^{(j,k)}(g)=\mathrm{tg}'^{(j,k)}(f)$, then the following chain of equalities holds
\allowdisplaybreaks
\begin{align*}
\beta^{\mathsf{C}}_{k}\left(
g\circ'^{(j)}f
\right)&=
\mathrm{sc}^{(k,n)}\left(
g\circ'^{(j)}f
\right)
\tag{1}
\\&=
\mathrm{sc}^{(k,n)}\left(
g\circ^{(j)}f
\right)
\tag{2}
\\&=
\mathrm{sc}^{(k,n)}\left(
g
\right)\circ^{(j)}\mathrm{sc}^{(k,n)}\left(
f
\right)
\tag{3}
\\&=
\beta^{\mathsf{C}}_{k}(g)\circ^{(j)}\beta^{\mathsf{C}}_{k}(f).
\tag{4}
\end{align*}

The first equality unravels the definition of the $k$-th component of $\beta^{\mathsf{C}}$; the second equality unravels the definition of the partial operation of binary composition $\circ'^{(j)}$ in $\mathrm{Gd}^{(n)}(\mathrm{Uf}^{(n)}(\mathsf{C}))$; the third equality follows from the fact that, by item~(MS3) in Definition~\ref{DnCatMS} it is the case $\mathrm{sc}^{(k,n)}(g\circ^{(j)}f)=\mathrm{sc}^{(k,n)}(g)\circ^{(j)}\mathrm{sc}^{(k,n)}(f)$; finally, the last equality recovers the $k$-th component of $\beta^{\mathsf{C}}$.

For the particular case at the $n$-th sort, for every $0\leq k<n$, the following chain of equalities holds.
\allowdisplaybreaks
\begin{align*}
\beta^{\mathsf{C}}_{n}\left(
g\circ'^{(j)}f
\right)&=
g\circ'^{(j)}f
\tag{1}
\\&=
g\circ^{(j)}f
\tag{2}
\\&=
\beta^{\mathsf{C}}_{k}(g)\circ^{(j)}\beta^{\mathsf{C}}_{k}(f).
\tag{3}
\end{align*}

The first equality unravels the definition of the $n$-th component of $\beta^{\mathsf{C}}$. Let us recall that 
$\beta^{\mathsf{C}}_{n}=\mathrm{Id}_{C_{n}}$; the second equality unravels the definition of the partial operation of binary composition $\circ'^{(j)}$ in $\mathrm{Gd}^{(n)}(\mathrm{Uf}^{(n)}(\mathsf{C}))$;  finally, the last equality recovers the $n$-th component of $\beta^{\mathsf{C}}$. Let us recall that 
$\beta^{\mathsf{C}}_{n}=\mathrm{Id}_{C_{n}}$.

All in all,  $\beta^{\mathsf{C}}$ is a many-sorted $n$-functor from  $\mathrm{Gd}^{(n)}(\mathrm{Uf}^{(n)}(\mathsf{C}))$ to $\mathsf{C}$,
$$\beta^{\mathsf{C}}\colon \mathrm{Gd}^{(n)}(\mathrm{Uf}^{(n)}(\mathsf{C}))\mor \mathsf{C}.$$

Claim~\ref{CGradUnifD} follows.

We next prove that, for every many-sorted $n$-category $\mathsf{C}$ in $\mathsf{nCat}^{\mathrm{MS}}$, the many-sorted $n$-functors $\alpha^{\mathsf{C}}$ and $\beta^{\mathsf{C}}$ are mutually inverse morphisms in $\mathsf{nCat}^{\mathrm{MS}}$.

\begin{claim}\label{CGradUnifE} For every many-sorted $n$-category $\mathsf{C}$ in $\mathsf{nCat}^{\mathrm{MS}}$, 
\allowdisplaybreaks
\begin{align*}
\beta^{\mathsf{C}}\circ\alpha^{\mathsf{C}}&=\mathrm{Id}_{\mathsf{C}},
&
\alpha^{\mathsf{C}}\circ\beta^{\mathsf{C}}&=\mathrm{Id}_{\mathrm{Gd}^{(n)}(\mathrm{Uf}^{(n)}(\mathsf{C}))}.
\end{align*}
\end{claim}
For every $0\leq k<n$, if $f\in C_{k}$, the following chain of equalities holds
\allowdisplaybreaks
\begin{align*}
\beta^{\mathsf{C}}_{k}\left(
\alpha^{\mathsf{C}}_{k}\left(
f
\right)\right)
&=
\mathrm{sc}^{(k,n)}\left(
\alpha^{\mathsf{C}}_{k}\left(
f
\right)\right)
\tag{1}
\\&=
\mathrm{sc}^{(k,n)}\left(
\mathrm{i}^{(n,k)}\left(
f
\right)\right)
\tag{2}
\\&=
f.
\tag{3}
\end{align*}

The first equality unravels the definition of the $k$-th component of $\beta^{\mathsf{C}}$; the second equality unravels the definition of the $k$-th component of $\alpha^{\mathsf{C}}$; finally, the last equality follows from the fact that, by item~(MS5) in Definition~\ref{DnCatMS}, $\mathrm{sc}^{(k,n)}(\mathrm{i}^{(n,k)}(f))=f$.

For the particular case at the $n$-th sort, if $f\in C_{n}$, then the following chain of equalities holds
\allowdisplaybreaks
\begin{align*}
\beta^{\mathsf{C}}_{n}\left(
\alpha^{\mathsf{C}}_{n}\left(
f
\right)\right)
&=
\alpha^{\mathsf{C}}_{n}(f)
\tag{1}
\\&=
f.
\tag{2}
\end{align*}

The first equality unravels the definition of the $n$-th component of $\beta^{\mathsf{C}}$. Let us recall that $\beta^{\mathsf{C}}_{n}=\mathrm{Id}_{C_{n}}$; finally, the second equality unravels the definition of the $n$-th component of $\alpha^{\mathsf{C}}$. Let us recall that $\alpha^{\mathsf{C}}_{n}=\mathrm{Id}_{C_{n}}$.

Thus, $\beta^{\mathsf{C}}\circ\alpha^{\mathsf{C}}=\mathrm{Id}_{\mathsf{C}}$.

Conversely, for every $0\leq k<n$, if $f\in C'_{k}$, then the following chain of equalities holds
\allowdisplaybreaks
\begin{align*}
\alpha^{\mathsf{C}}_{k}\left(
\beta^{\mathsf{C}}_{k}\left(
f
\right)\right)
&=
\mathrm{i}^{(n,k)}\left(
\beta^{\mathsf{C}}_{k}\left(
f
\right)\right)
\tag{1}
\\&=
\mathrm{i}^{(k,n)}\left(
\mathrm{sc}^{(n,k)}\left(
f
\right)\right)
\tag{2}
\\&=
f.
\tag{3}
\end{align*}

The first equality unravels the definition of the $k$-th component of $\alpha^{\mathsf{C}}$; the second equality unravels the definition of the $k$-th component of $\alpha^{\mathsf{C}}$; finally, the last equality follows from the fact that $f\in C'_{k}$, hence $\mathrm{i}^{(k,n)}(\mathrm{sc}^{(n,k)}(f))=f$.

For the particular case at the $n$-th sort, if $f\in C_{n}$, then the following chain of equalities holds
\allowdisplaybreaks
\begin{align*}
\alpha^{\mathsf{C}}_{n}\left(
\beta^{\mathsf{C}}_{n}\left(
f
\right)\right)
&=
\beta^{\mathsf{C}}_{n}(f)
\tag{1}
\\&=
f.
\tag{2}
\end{align*}

The first equality unravels the definition of the $n$-th component of $\alpha^{\mathsf{C}}$. Let us recall that $\alpha^{\mathsf{C}}_{n}=\mathrm{Id}_{C_{n}}$; finally, the second equality unravels the definition of the $n$-th component of $\beta^{\mathsf{C}}$. Let us recall that $\beta^{\mathsf{C}}_{n}=\mathrm{Id}_{C_{n}}$.

Thus, $\alpha^{\mathsf{C}}\circ\beta^{\mathsf{C}}=\mathrm{Id}_{\mathrm{Gd}^{(n)}(\mathrm{Uf}^{(n)}(\mathsf{C}))}$.

Claim~\ref{CGradUnifE} follows.

We next prove that the arrows $\alpha^{\mathsf{C}}$, when $\mathsf{C}$ ranges over the objects of $\mathsf{nCat}^{\mathrm{MS}}$, are the components of a natural transformation.

\begin{claim}\label{CGradUnifF}
For every many-sorted $n$-functor $F\colon \mathsf{C}\mor \mathsf{D}$, the square of many-sorted $n$-functors in Figure~\ref{FAlphaFun} commutes.
\end{claim}
\begin{figure}
\begin{center}
\begin{tikzpicture}
[ACliment/.style={-{To [angle'=45, length=5.75pt, width=4pt, round]}}
, scale=1, 
AClimentD/.style={dogble eqgal sign distance,
-implies
}
]

\node[] (a) at (0,0) [] {$\mathsf{C}$};
\node[] (b) at (3,0) [] {$\mathrm{Gd}^{(n)}(\mathrm{Uf}^{(n)}(\mathsf{C}))$};
\node[] (c) at (0,-1.5) [] {$\mathsf{D}$};
\node[] (d) at (3,-1.5) [] {$\mathrm{Gd}^{(n)}(\mathrm{Uf}^{(n)}(\mathsf{D}))$};

\draw[ACliment] (a) to node [above] {$\scriptstyle \alpha^{\mathsf{C}}$} (b);
\draw[ACliment] (a) to node [left] {$\scriptstyle F$} (c);
\draw[ACliment] (c) to node [below] {$\scriptstyle \alpha^{\mathsf{D}}$} (d);
\draw[ACliment] (b) to node [right] {$\scriptstyle \mathrm{Gd}^{(n)}(\mathrm{Uf}^{(n)}(F))$} (d);
\end{tikzpicture}
\end{center}
\caption{The natural transformation $\alpha$ on many-sorted $n$-functors.}
\label{FAlphaFun}
\end{figure}

Let $\mathsf{C}$ and $\mathsf{D}$ be the $n$-categories under consideration, where
\allowdisplaybreaks
\begin{align*}
\mathsf{C}&=
\left(\left(C_{k}
\right)_{k\in n+1}, 
\left(
\circ^{(j)\mathsf{C}},\mathrm{sc}^{(j,k)\mathsf{C}},\mathrm{tg}^{(j,k)\mathsf{C}},\mathrm{i}^{(k,j)\mathsf{C}}
\right)_{(j,k)\in\coprod_{k\in n+1}k}
\right)\,\text{ and}\\
\mathsf{D}&=
\left(\left(
D_{k}
\right)_{k\in n+1}, 
\left(
\circ^{(j)\mathsf{D}},\mathrm{sc}^{(j,k)\mathsf{D}},\mathrm{tg}^{(j,k)\mathsf{D}},\mathrm{i}^{(k,j)\mathsf{D}}
\right)_{(j,k)\in\coprod_{k\in n+1}k}
\right),
\end{align*}
and let 
\allowdisplaybreaks
\begin{align*}
\mathrm{Gd}^{(n)}\left(
\mathrm{Uf}^{(n)}\left(
\mathsf{C}
\right)\right)&=
\left(\left(
C'_{k}
\right)_{k\in n+1}, 
\left(
\circ'^{(j)\mathsf{C}},\mathrm{sc}'^{(j,k)\mathsf{C}},\mathrm{tg}'^{(j,k)\mathsf{C}},\mathrm{i}'^{(k,j)\mathsf{C}}
\right)_{(j,k)\in\coprod_{k\in n+1}k}\right)\,\text{ and}\\
\mathrm{Gd}^{(n)}\left(\mathrm{Uf}^{(n)}\left(
\mathsf{D}
\right)\right)&=
\left(\left(
D'_{k}
\right)_{k\in n+1},
\left(
\circ'^{(j)\mathsf{D}},\mathrm{sc}'^{(j,k)\mathsf{D}}, \mathrm{tg}'^{(j,k)\mathsf{D}},\mathrm{i}'^{(k,j)\mathsf{D}}
\right)_{(j,k)\in\coprod_{k\in n+1}k}\right)
\end{align*}
be the corresponding graduations of the unifications of $\mathsf{C}$ and $\mathsf{D}$, respectively.

Let $F=(F_{k})_{k\in n+1}$ be a many-sorted $n$-functor from $\mathsf{C}$ to $\mathsf{D}$. Let us recall that, according to Definition~\ref{DUnifFun}, its unification is given by $\mathrm{Uf}^{(n)}(F)=F_{n}$. Moreover, according to Definition~\ref{DGradFun}, the graduation of the unification of $F$ is given by $\mathrm{Gd}^{(n)}(\mathrm{Uf}^{(n)}(F))=(\mathrm{Gd}^{(n)}(\mathrm{Uf}^{(n)}(F))_{k})_{k\in n+1}$ where, for every $k\in n+1$, 
$$
\mathrm{Gd}^{(n)}(\mathrm{Uf}^{(n)}(F))_{k}=
\begin{cases}
F_{n}\bigr|_{C'_{k}}^{D'_{k}},&\mbox{if }k\in n;
\\
F_{n},&\mbox{if } k=n.
\end{cases}
$$

Let $\alpha^{\mathsf{C}}=(\alpha^{\mathsf{C}}_{k})_{k\in n+1}$ be a many-sorted $n$-functor from $\mathsf{C}$ to $\mathrm{Gd}^{(n)}(\mathrm{Uf}^{(n)}(\mathsf{C}))$ and $\alpha^{\mathsf{D}}=(\alpha^{\mathsf{D}}_{k})_{k\in n+1}$ a many-sorted $n$-functor from $\mathsf{D}$ to $\mathrm{Gd}^{(n)}(\mathrm{Uf}^{(n)}(\mathsf{D}))$. Let us recall that, for every $k\in n+1$, $\alpha^{\mathsf{C}}_{k}$ and $\alpha^{\mathsf{D}}_{k}$ are defined as 
$$
\alpha^{\mathsf{C}}_{k}=
\begin{cases}
\mathrm{i}^{(n,k)\mathsf{C}},
&\mbox{if }k\in n;
\\
\mathrm{Id}_{C_{n}},
&\mbox{if } k=n;
\end{cases}
\qquad\qquad
\alpha^{\mathsf{D}}_{k}=
\begin{cases}
\mathrm{i}^{(n,k)\mathsf{D}},
&\mbox{if }k\in n;
\\
\mathrm{Id}_{D_{n}},
&\mbox{if } k=n.
\end{cases}
$$

We represent in the  diagram of Figure~\ref{FAlphaNTnCat} all the different sets and mappings under consideration.

\begin{figure}
\begin{center}
\begin{tikzpicture}
[ACliment/.style={-{To [angle'=45, length=5.75pt, width=4pt, round]}}]
    
    \node (C0) at (0,0) {$C_{0}$};
    \node (C1) at (1,1) {$C_{1}$};
    \node (C2) at (2,2) {$\cdots$};
    \node (C3) at (3,3) {$C_{n}$};  
    \node (D0) at (0,-2.5) {$D_{0}$};
    \node (D1) at (1,-1.5) {$D_{1}$};
    \node (D2) at (2,-.5) {$\cdots$};
    \node (D3) at (3,.5) {$D_{n}$};
    \node (C0p) at (6,0) {$C'_{0}$};
    \node (C1p) at (7,1) {$C'_{1}$};
    \node (C2p) at (8,2) {$\cdots$};
    \node (C3p) at (9,3) {$C'_{n}$};         
    \node (D0p) at (6,-2.5) {$D'_{0}$};
    \node (D1p) at (7,-1.5) {$D'_{1}$};
    \node (D2p) at (8,-.5) {$\cdots$};
    \node (D3p) at (9,.5) {$D'_{n}$};
    
    \begin{scope}[transparency group, opacity=0.2]
    \draw[ACliment] (C1.205) -- (C0.65);
    \draw[ACliment] (C0.45) -- (C1.225);
    \draw[ACliment] (C1.245) -- (C0.385);
    
    \draw[ACliment] (C2.205) -- (C1.65);
    \draw[ACliment] (C1.45) -- (C2.225);
    \draw[ACliment] (C2.245) -- (C1.385);
    
     \draw[ACliment] (C3.205) -- (C2.65);
    \draw[ACliment] (C2.45) -- (C3.225);
    \draw[ACliment] (C3.245) -- (C2.385);
    \end{scope}
    
    \begin{scope}[transparency group, opacity=0.2]
    \draw[ACliment] (D1.205) -- (D0.65);
    \draw[ACliment] (D0.45) -- (D1.225);
    \draw[ACliment] (D1.245) -- (D0.385);
    
    \draw[ACliment] (D2.205) -- (D1.65);
    \draw[ACliment] (D1.45) -- (D2.225);
    \draw[ACliment] (D2.245) -- (D1.385);
    
     \draw[ACliment] (D3.205) -- (D2.65);
    \draw[ACliment] (D2.45) -- (D3.225);
    \draw[ACliment] (D3.245) -- (D2.385);
    \end{scope}
    
    \begin{scope}[transparency group, opacity=0.2]
    \draw[ACliment] (C1p.205) -- (C0p.65);
    \draw[ACliment] (C0p.45) -- (C1p.225);
    \draw[ACliment] (C1p.245) -- (C0p.385);
    
    \draw[ACliment] (C2p.205) -- (C1p.65);
    \draw[ACliment] (C1p.45) -- (C2p.225);
    \draw[ACliment] (C2p.245) -- (C1p.385);
    
     \draw[ACliment] (C3p.205) -- (C2p.65);
    \draw[ACliment] (C2p.45) -- (C3p.225);
    \draw[ACliment] (C3p.245) -- (C2p.385);
    \end{scope}
    
     \begin{scope}[transparency group, opacity=0.2]
    \draw[ACliment] (D1p.205) -- (D0p.65);
    \draw[ACliment] (D0p.45) -- (D1p.225);
    \draw[ACliment] (D1p.245) -- (D0p.385);
    
    \draw[ACliment] (D2p.205) -- (D1p.65);
    \draw[ACliment] (D1p.45) -- (D2p.225);
    \draw[ACliment] (D2p.245) -- (D1p.385);
    
     \draw[ACliment] (D3p.205) -- (D2p.65);
    \draw[ACliment] (D2p.45) -- (D3p.225);
    \draw[ACliment] (D3p.245) -- (D2p.385);
    \end{scope}
    
    \draw[ACliment] (C0) to node [midway, fill=white] {$\scriptstyle \alpha^{\mathsf{C}}_{0}$} (C0p);
    \draw[ACliment] (C1) to node [midway, fill=white] {$\scriptstyle\alpha^{\mathsf{C}}_{1}$} (C1p);
    \draw[ACliment, white, text=black] (C2) to node [midway, fill=white] {$\scriptstyle\cdots$} (C2p);
    \draw[ACliment] (C3) to node [midway, fill=white] {$\scriptstyle \alpha^{\mathsf{C}}_{n}$} (C3p);
    
    \draw[ACliment] (D0) to node [midway, fill=white] {$\scriptstyle \alpha^{\mathsf{D}}_{0}$} (D0p);
    \draw[ACliment] (D1) to node [midway, fill=white] {$\scriptstyle\alpha^{\mathsf{D}}_{1}$} (D1p);
    \draw[ACliment, white, text=black] (D2) to node [midway, fill=white] {$\scriptstyle\cdots$} (D2p);
    \draw[ACliment] (D3) to node [midway, fill=white] {$\scriptstyle \alpha^{\mathsf{D}}_{n}$} (D3p);    
    
    \draw[ACliment] (C0) to node [midway, fill=white] {$\scriptstyle F_{0}$} (D0);
    \draw[ACliment] (C1) to node [midway, fill=white] {$\scriptstyle F_{1}$} (D1);
    \draw[ACliment, white, text=black] (C2) to node [midway, fill=white] {$\scriptstyle\cdots$} (D2);
    \draw[ACliment] (C3) to node [midway, fill=white] {$\scriptstyle F_{n}$} (D3);
    
    \draw[ACliment] (C0p) to node [midway, fill=white] {$\scriptstyle \mathrm{Gd}^{(n)}(\mathrm{Uf}^{(n)}(F))_{0}$} (D0p);
    \draw[ACliment] (C1p) to node [midway, fill=white] {$\scriptstyle \mathrm{Gd}^{(n)}(\mathrm{Uf}^{(n)}(F))_{1}$} (D1p);
    \draw[ACliment, white, text=black] (C2p) to node [midway, fill=white] {$\scriptstyle\cdots$} (D2p);
    \draw[ACliment] (C3p) to node [midway, fill=white] {$\scriptstyle \mathrm{Gd}^{(n)}(\mathrm{Uf}^{(n)}(F))_{n}$} (D3p);
       
    \node (C0) at (0,0) {$C_{0}$};
    \node (C1) at (1,1) {$C_{1}$};
    \node (C2) at (2,2) {$\cdots$};
    \node (C3) at (3,3) {$C_{n}$};  
    \node (D0) at (0,-2.5) {$D_{0}$};
    \node (D1) at (1,-1.5) {$D_{1}$};
    \node (D2) at (2,-.5) {$\cdots$};
    \node (D3) at (3,.5) {$D_{n}$};
    \node (C0p) at (6,0) {$C'_{0}$};
    \node (C1p) at (7,1) {$C'_{1}$};
    \node (C2p) at (8,2) {$\cdots$};
    \node (C3p) at (9,3) {$C'_{n}$};         
    \node (D0p) at (6,-2.5) {$D'_{0}$};
    \node (D1p) at (7,-1.5) {$D'_{1}$};
    \node (D2p) at (8,-.5) {$\cdots$};
    \node (D3p) at (9,.5) {$D'_{n}$}; 
\end{tikzpicture} 
\end{center}
\caption{The natural transformation $\alpha$ on many-sorted $n$-categories.}
\label{FAlphaNTnCat}
\end{figure}

In order to prove Claim~\ref{CGradUnifE} we need to check that, for every $k\in n+1$, the diagram in Figure~\ref{FAlphaNTk} commutes.

\begin{figure}
\begin{center}
\begin{tikzpicture}
[ACliment/.style={-{To [angle'=45, length=5.75pt, width=4pt, round]}}
, scale=1, 
AClimentD/.style={dogble eqgal sign distance,
-implies
}
]

\node[] (a) at (0,0) [] {$C_{k}$};
\node[] (b) at (3,0) [] {$\mathrm{Gd}^{(n)}(\mathrm{Uf}^{(n)}(\mathsf{C}))_{k}$};
\node[] (c) at (0,-1.5) [] {$D_{k}$};
\node[] (d) at (3,-1.5) [] {$\mathrm{Gd}^{(n)}(\mathrm{Uf}^{(n)}(\mathsf{D}))_{k}$};

\draw[ACliment] (a) to node [above] {$\scriptstyle \alpha^{\mathsf{C}}_{k}$} (b);
\draw[ACliment] (a) to node [left] {$\scriptstyle F_{k}$} (c);
\draw[ACliment] (c) to node [below] {$\scriptstyle \alpha^{\mathsf{D}}_{k}$} (d);
\draw[ACliment] (b) to node [right] {$\scriptstyle \mathrm{Gd}^{(n)}(\mathrm{Uf}^{(n)}(F))_{k}$} (d);
\end{tikzpicture}
\end{center}
\caption{The natural transformation $\alpha$ at layer $k$ of a many-sorted $n$-category.}
\label{FAlphaNTk}
\end{figure}

Let $k\in n$ and $f\in C_{k}$, then the following chain of equalities holds
\allowdisplaybreaks
\begin{align*}
\left(
\mathrm{Gd}^{(n)}\circ\mathrm{Uf}^{(n)}(F)
\right)_{k}\left(
\alpha^{\mathsf{C}}_{k}(f)
\right)
&=
F_{n}\bigr|_{C'_{k}}^{D'_{k}}\left(
\alpha^{\mathsf{C}}_{k}(f)
\right)
\tag{1}
\\&=
F_{n}\bigr|_{C'_{k}}^{D'_{k}}\left(
\mathrm{i}^{(n,k)\mathsf{C}}\left(
f
\right)\right)
\tag{2}
\\&=
F_{n}\left(
\mathrm{i}^{(n,k)\mathsf{C}}\left(
f
\right)\right)
\tag{3}
\\&=
\mathrm{i}^{(n,k)\mathsf{D}}\left(
F_{k}\left(
f
\right)\right)
\tag{4}
\\&=
\alpha^{\mathsf{D}}_{k}\left(
F_{k}\left(
f
\right)\right).
\tag{5}
\end{align*}

The first equality unravels the $k$-th component of the many-sorted $n$-functor $\mathrm{Gd}^{(n)}\circ\mathrm{Uf}^{(n)}$;
the second equality unravels the $k$-th component of the many-sorted $n$-functor $\alpha^{\mathsf{C}}$; the third equality applies the definition of a restriction of a mapping; the fourth equation follows from the fact that, since $F$ is a many-sorted $n$-functor, by item~(Msii) in Definition~\ref{DnCatMS}, $F_{n}(\mathrm{i}^{(n,k)\mathsf{C}}(f))=\mathrm{i}^{(n,k)\mathsf{D}}(F_{k}(f))$; finally, the last equality recovers the $k$-th component of the many-sorted $n$-functor $\alpha^{\mathsf{D}}$.

Finally, for $k=n$ and $f\in C_{n}$, the following chain of equalities holds
\allowdisplaybreaks
\begin{align*}
\left(
\mathrm{Gd}^{(n)}\circ\mathrm{Uf}^{(n)}
\right)_{n}\left(
\alpha^{\mathsf{C}}_{n}\left(
f
\right)\right)
&=
F_{n}\left(
\alpha^{\mathsf{C}}_{n}\left(
f
\right)\right)
\tag{1}
\\&=
F_{n}\left(
\mathrm{Id}_{C_{n}}\left(
f
\right)\right)
\tag{2}
\\&=
F_{n}(f)
\tag{3}
\\&=
\mathrm{Id}_{D_{n}}\left(
F_{n}\left(
f
\right)\right)
\tag{4}
\\&=
\alpha^{\mathsf{D}}_{n}\left(
F_{n}\left(
f
\right)\right).
\tag{5}
\end{align*}

The first equality unravels the $n$-th component of the many-sorted $n$-functor $\mathrm{Gd}^{(n)}\circ\mathrm{Uf}^{(n)}$;
the second equality unravels the $n$-th component of the many-sorted $n$-functor $\alpha^{\mathsf{C}}$; 
the third equality applies the definition of an identity mapping; the fourth equality applies the definition of an identity mapping; finally, the last equality recovers the $n$-th component of the many-sorted $n$-functor $\alpha^{\mathsf{D}}$.

Claim~\ref{CGradUnifF} follows.

So, considering the foregoing, we can affirm that $\alpha$ is a natural transformation from $\mathrm{Id}_{\mathsf{nCat}^{\mathrm{MS}}}$ to  $\mathrm{Gd}^{(n)}\circ\mathrm{Uf}^{(n)}$.

\begin{figure}
\begin{center}
\begin{tikzpicture}
[ACliment/.style={-{To [angle'=45, length=5.75pt, width=4pt, round]}}
, scale=0.8, 
AClimentD/.style={double equal sign distance,
-implies
}
]
\node[] (u) at (0,0) [] {$\mathsf{nCat}^{\mathrm{MS}}$};
\node[] (v) at (4,0) [] {$\mathsf{nCat}^{\mathrm{MS}}$};
\draw[ACliment, bend left] (u) to node [above]	{$\mathrm{Id}_{\mathsf{nCat}^{\mathrm{MS}}}$} (v);
\draw[ACliment, bend right] (u) to node [below] {$\mathrm{Gd}^{(n)}\circ\mathrm{Uf}^{(n)}$} (v);
\node[] (a) at (2,.5) [] {};
\node[] (b) at (2,-.5) [] {};
\draw[AClimentD]  (a) to node [right]{$\alpha$} (b);
\end{tikzpicture}
\end{center}
\caption{The natural transformation $\alpha$.}
\label{FAlphaNT}
\end{figure}

Let us note that each many-sorted $n$-category $\mathsf{C}$ in $\mathsf{nCat}^{\mathrm{MS}}$ determines, by Claim~\ref{CGradUnifC}, the many-sorted $n$-functor $\alpha^{\mathsf{C}}$ which, by Claim~\ref{CGradUnifE}, is an isomorphism. Thus, $\alpha\colon \mathrm{Id}_{\mathsf{nCat}^{\mathrm{MS}}}\cel\mathrm{Gd}^{(n)}\circ\mathrm{Uf}^{(n)}$ is a natural isomorphism.
\end{proof}

From the above, Proposition~\ref{Eqvt} and Definition~\ref{lali} we obtain the following corollary.

\begin{corollary}\label{CDnCatEquiv} 
The functor $\mathrm{Uf}^{(n)}$ is a left adjoint left inverse of the functor $\mathrm{Gd}^{(n)}$. Hence the categories $\mathsf{nCat}$ and $\mathsf{nCat}^{\mathrm{MS}}$ are equivalent.
\end{corollary}



\section{
\texorpdfstring
{From  single-sorted to  many-sorted  \(\omega\)-categories}
{From  single-sorted to  many-sorted  Omega-categories}
}

In this section we associate to a single-sorted $\omega$-category a many-sorted $\omega$-category by means of the notion of \(n\)-graduation and the universal property of \(\w\mathsf{Cat}^{\mathrm{MS}}\). Moreover, we also associate to a single-sorted $\omega$-functor between single-sorted $\omega$-categories a many-sorted $\omega$-functor between the associated many-sorted $\omega$-categories. We end this section by explicitly giving the definition of the \(\omega\)-graduation functor.

\begin{proposition}\label{PGradSq}
Let \(m,n\in\omega\) with \(m<n\). The square
\[
\xymatrix{
\mathsf{mCat}
  \ar[d]_{\mathrm{Gd}^{(m)}}
  &
\mathsf{nCat}
  \ar[l]_{U^{(m,n)}}
  \ar[d]^{\mathrm{Gd}^{(n)}}
\\
\mathsf{mCat}^{\mathrm{MS}}
  &
\mathsf{nCat}^{\mathrm{MS}}
  \ar[l]^{U^{(m,n)}_{\mathrm{MS}}}
}
\]
commutes.
\end{proposition}

\begin{proof}
Let \(\mathsf{C}=(C,(\xi_{k})_{k\in n})\) be a single-sorted \(n\)-category, with \((\xi_{k})_{k\in n}=(\#^{k},\mathrm{sc}^{k},\mathrm{tg}^{k})_{k\in n}\).
We want to prove that
\[
\mathrm{Gd}^{(m)}(\mathsf{C}^{<m})=\mathrm{Gd}^{(n)}(\mathsf{C})^{<m}.
\]

Following Definition~\ref{DUnderlying}, the \(m\)-category \(\mathsf{C}^{<m}\) is \((C_{m},(\xi^{<m}_{k})_{k\in m})\) with \((\xi^{<m}_{k})_{k\in m}=(\#^{k}\bigr|^{C_{m}}_{C_{m}\times C_{m}},\mathrm{sc}^{k}\bigr|^{C_{m}}_{C_{m}},\mathrm{tg}^{k}\bigr|^{C_{m}}_{C_{m}})\). Now, from Definition~\ref{DGrad}, its graduation is the many-sorted \(m\)-category \((((C_{m})_{k})_{k\in m+1}, \zeta)\) where \(\zeta=(\zeta_{j,k})_{(j,k)\in\coprod_{k\in m+1}k}\) with
\[
\mathrm{Gd}^{(m)}(\mathsf{C}^{<m})
=
\left(
\left(\left(C_{m}\right)_{k}\right)_{k\in m+1},
\left(
\circ^{(j)}, \mathrm{sc}^{(j,k)}, \mathrm{tg}^{(j,k)}, \mathrm{i}^{(k,j)}
\right)_{(j,k)\in\coprod_{k\in m+1}k}
\right)
\]
Moreover, for every \(0\leq j<k\leq m\), the structural operations of the graduation are given by
\begin{align*}
\circ^{(j)}&=\#^{j}\bigr|^{C_{m}}_{C_{m}\times C_{m}}\bigr|^{(C_{m})_{k}}_{(C_{m})_{k}\times(C_{m})_{k}};
&\mathrm{i}^{(k,j)}&=\mathrm{in}^{(C_{m})_{j},(C_{m})_{k}};
\\
\mathrm{sc}^{(j,k)}&=\mathrm{sc}^{j}\bigr|^{C_{m}}_{C_{m}}\bigr|^{(C_{m})_{j}}_{(C_{m})_{k}};
&\mathrm{tg}^{(j,k)}&=\mathrm{tg}^{j}\bigr|^{C_{m}}_{C_{m}}\bigr|^{(C_{m})_{j}}_{(C_{m})_{k}}.
\end{align*}

All in all, the following diagram represents the many-sorted \(m\)-category \(\mathrm{Gd}^{(m)}(\mathsf{C}^{<m})\).
\[
\xymatrixcolsep={18ex}
\xymatrixrowsep={12ex}
\scalebox{0.75}{\xymatrix{
(C_{m})_{0}
  \ar[r]|-{\mathrm{in}^{(C_{m})_{0},(C_{m})_{1}}}
  &
(C_{m})_{1}
  \ar@<-3ex>[l]_{\mathrm{sc}^{0}\bigr|^{C_{m}}_{C_{m}}\bigr|^{(C_{m})_{0}}_{(C_{m})_{1}}}
  \ar[r]|-{\mathrm{in}^{(C_{m})_{1},(C_{m})_{2}}}
  \ar@<3ex>[l]^{\mathrm{tg}^{0}\bigr|^{C_{m}}_{C_{m}}\bigr|^{(C_{m})_{0}}_{(C_{m})_{1}}}
  &
(C_{m})_{2}
  \ar@<-3ex>[l]_{\mathrm{sc}^{1}\bigr|^{C_{m}}_{C_{m}}\bigr|^{(C_{m})_{1}}_{(C_{m})_{2}}}
  \ar[r]|-{\mathrm{in}^{(C_{m})_{2},(C_{m})_{3}}}
  \ar@<3ex>[l]^{\mathrm{tg}^{1}\bigr|^{C_{m}}_{C_{m}}\bigr|^{(C_{m})_{1}}_{(C_{m})_{2}}}
  &
\cdots
  \ar@<-3ex>[l]_{\mathrm{sc}^{2}\bigr|^{C_{m}}_{C_{m}}\bigr|^{(C_{m})_{2}}_{(C_{m})_{3}}}
  \ar[r]|-{\mathrm{in}^{(C_{m})_{m-1},(C_{m})_{m}}}
  \ar@<3ex>[l]^{\mathrm{tg}^{2}\bigr|^{C_{m}}_{C_{m}}\bigr|^{(C_{m})_{2}}_{(C_{m})_{3}}}
  &
(C_{m})_{m}
  \ar@<-3ex>[l]_{\mathrm{sc}^{m-1}\bigr|^{C_{m}}_{C_{m}}\bigr|^{(C_{m})_{m-1}}_{(C_{m})_{m}}}
  \ar@<3ex>[l]^{\mathrm{tg}^{m-1}\bigr|^{C_{m}}_{C_{m}}\bigr|^{(C_{m})_{m-1}}_{(C_{m})_{m}}}
}}
\]

Following Definition~\ref{DGrad}, the \(n\)-category \(\mathrm{Gd}^{(n)}(\mathsf{C})\) is
\[
\mathrm{Gd}^{(n)}(\mathsf{C})
=
\left(
\left(C_{k}\right)_{k\in n+1},
\left( \circ'^{(j)}, \mathrm{sc}'^{(j,k)}, \mathrm{tg}'^{(j,k)}, \mathrm{i}'^{(k,j)} \right)_{(j,k)\in\coprod_{k\in n+1}k}
\right)
\]
where, for every \(0\leq j<k\leq n\), the structural operations of the graduation are given by
\begin{align*}
\circ'^{(j)}&=\#^{j}\bigr|^{C_{k}}_{C_{k}\times C_{k}};
&\mathrm{i}'^{(k,j)}&=\mathrm{in}^{C_{j},C_{k}};
\\
\mathrm{sc}'^{(j,k)}&=\mathrm{sc}^{j}\bigr|^{C_{j}}_{C_{k}};
&\mathrm{tg}'^{(j,k)}&=\mathrm{tg}^{j}\bigr|^{C_{j}}_{C_{k}}.
\end{align*}
Now, from Definition~\ref{DUnderlyingMS}, its underlying many-sorted \(m\)-category is \(((C_{k})_{k\in m+1}, \zeta^{<m})\) with \(\zeta^{<m}=(\zeta_{j,k})_{(j,k)\in\coprod_{k\in m+1}k}\).

All in all, the following diagram represents the many-sorted \(m\)-category \(\mathrm{Gd}^{(m)}(\mathsf{C})^{<m}\).
\[
\xymatrixcolsep={18ex}
\xymatrixrowsep={12ex}
\scalebox{0.85}{\xymatrix{
C_{0}
  \ar[r]|-{\mathrm{in}^{C_{0},C_{1}}}
  &
C_{1}
  \ar@<-3ex>[l]_{\mathrm{sc}^{0}\bigr|^{C_{0}}_{C_{1}}}
  \ar[r]|-{\mathrm{in}^{C_{1},C_{2}}}
  \ar@<3ex>[l]^{\mathrm{tg}^{0}\bigr|^{C_{0}}_{C_{1}}}
  &
C_{2}
  \ar@<-3ex>[l]_{\mathrm{sc}^{1}\bigr|^{C_{1}}_{C_{2}}}
  \ar[r]|-{\mathrm{in}^{C_{2},C_{3}}}
  \ar@<3ex>[l]^{\mathrm{tg}^{1}\bigr|^{C_{1}}_{C_{2}}}
  &
\cdots
  \ar@<-3ex>[l]_{\mathrm{sc}^{2}\bigr|^{C_{2}}_{C_{3}}}
  \ar[r]|-{\mathrm{in}^{C_{m-1},C_{m}}}
  \ar@<3ex>[l]^{\mathrm{tg}^{2}\bigr|^{C_{2}}_{C_{3}}}
  &
C_{m}
  \ar@<-3ex>[l]_{\mathrm{sc}^{m-1}\bigr|^{C_{m-1}}_{C_{m}}}
  \ar@<3ex>[l]^{\mathrm{tg}^{m-1}\bigr|^{C_{m-1}}_{C_{m}}}
}}
\]

Therefore, in order to prove that \(\mathrm{Gd}^{(m)}(\mathsf{C}^{<m})=\mathrm{Gd}^{(n)}(\mathsf{C})^{<m}\), it suffices to note that, for every \(0\leq k\leq m\) \((C_{m})_{k}=C_{k}\) since,
\begin{align*}
\circ^{(j)}&=\#^{j}\bigr|^{C_{m}}_{C_{m}\times C_{m}}\bigr|^{(C_{m})_{k}}_{(C_{m})_{k}\times(C_{m})_{k}}=\#^{j}\bigr|^{(C_{m})_{k}}_{(C_{m})_{k}\times(C_{m})_{k}}=\circ'^{(j)};
\\
\mathrm{sc}^{(j,k)}&=\mathrm{sc}^{j}\bigr|^{C_{m}}_{C_{m}}\bigr|^{(C_{m})_{j}}_{(C_{m})_{k}}=\mathrm{sc}^{j}\bigr|^{(C_{m})_{j}}_{(C_{m})_{k}}=\mathrm{sc}'^{(j,k)};\mbox{ and}
\\
\mathrm{tg}^{(j,k)}&=\mathrm{tg}^{j}\bigr|^{C_{m}}_{C_{m}}\bigr|^{(C_{m})_{j}}_{(C_{m})_{k}}=\mathrm{tg}^{j}\bigr|^{(C_{m})_{j}}_{(C_{m})_{k}}=\mathrm{tg}'^{(j,k)}.
\end{align*}

Now, let \(F\) be a single-sorted \(n\)-functor from \(\mathsf{C}\) to \(\mathsf{C}'\) where \(\mathsf{C}=(C,(\xi_{k})_{k\in n})\) with \((\xi_{k})_{k\in n}=(\#^{k},\mathrm{sc}^{k},\mathrm{tg}^{k})_{k\in n}\) and \(\mathsf{C}'=(C',(\xi'_{k})_{k\in n})\) with \((\xi'_{k})_{k\in n}=(\#'^{k},\mathrm{sc}'^{k},\mathrm{tg}'^{k})_{k\in n}\).

We want to prove that
\[
\mathrm{Gd}^{(m)}(F^{<m})=\mathrm{Gd}^{(n)}(F)^{<m}.
\]

Following Definition~\ref{DUnderlying}, the single-sorted \(m\)-functor \(F^{<m}\) is \(F\bigr|^{C'_{m}}_{C_{m}}\). Now, from Definition~\ref{DGrad}, its graduation is the many-sorted \(m\)-functor 
\[
\left(
F\bigr|^{C'_{m}}_{C_{m}}\bigr|^{C'_{k}}_{C_{k}}
\right)_{k\in m+1}.
\]

All in all, the following diagram represents the many-sorted \(m\)-functor \(\mathrm{Gd}^{(m)}(F^{<m})\).

\[
\xymatrixcolsep={12ex}
\scalebox{0.8}{\xymatrix{
(C_{m})_{0}
  \ar[d]_{F\bigr|^{C'_{m}}_{C_{m}}\bigr|^{(C'_{m})_{0}}_{(C_{m})_{0}}}
  \ar[r]
  &
(C_{m})_{1}
  \ar[d]^{F\bigr|^{C'_{m}}_{C_{m}}\bigr|^{(C'_{m})_{1}}_{(C_{m})_{1}}}
  \ar@<-1ex>[l]
  \ar[r]
  \ar@<1ex>[l]
  &
(C_{m})_{2}
  \ar[d]^{F\bigr|^{C'_{m}}_{C_{m}}\bigr|^{(C'_{m})_{2}}_{(C_{m})_{2}}}
  \ar@<-1ex>[l]
  \ar[r]
  \ar@<1ex>[l]
  &
\cdots
  \ar@<-1ex>[l]
  \ar[r]
  \ar@<1ex>[l]
  &
(C_{m})_{n}
  \ar[d]^{F\bigr|^{C'_{m}}_{C_{m}}\bigr|^{(C'_{m})_{m}}_{(C_{m})_{m}}}
  \ar@<-1ex>[l]
  \ar@<1ex>[l]
\\
(C'_{m})_{0}
  \ar[r]
  &
(C'_{m})_{1}
  \ar@<-1ex>[l]
  \ar[r]
  \ar@<1ex>[l]
  &
(C'_{m})_{2}
  \ar@<-1ex>[l]
  \ar[r]
  \ar@<1ex>[l]
  &
\cdots
  \ar@<-1ex>[l]
  \ar[r]
  \ar@<1ex>[l]
  &
(C'_{m})_{n}
  \ar@<-1ex>[l]
  \ar@<1ex>[l]
}}
\]

Following Definition~\ref{DGrad}, the many-sorted \(n\)-functor \(\mathrm{Gd}^{(n)}(F)\) is \((F\bigr|^{C'_{k}}_{C_{k}})_{k\in n+1}\). Now, from Definition~\ref{DUnderlyingMS}, its underlying many-sorted \(m\)-functor is \(\left(F\bigr|^{C'_{k}}_{C_{k}}\right)_{k\in m+1}\).

All in all, the following diagram represents the many-sorted \(m\)-functor \(\mathrm{Gd}^{(m)}(F)^{<m}\).

\[
\xymatrixcolsep={12ex}
\scalebox{0.95}{\xymatrix{
C_{0}
  \ar[d]_{F\bigr|^{C'_{0}}_{C_{0}}}
  \ar[r]
  &
C_{1}
  \ar[d]^{F\bigr|^{C'_{1}}_{C_{1}}}
  \ar@<-1ex>[l]
  \ar[r]
  \ar@<1ex>[l]
  &
C_{2}
  \ar[d]^{F\bigr|^{C'_{2}}_{C_{2}}}
  \ar@<-1ex>[l]
  \ar[r]
  \ar@<1ex>[l]
  &
\cdots
  \ar@<-1ex>[l]
  \ar[r]
  \ar@<1ex>[l]
  &
C_{n}
  \ar[d]^{F\bigr|^{C'_{m}}_{C_{m}}}
  \ar@<-1ex>[l]
  \ar@<1ex>[l]
\\
C'_{0}
  \ar[r]
  &
C'_{1}
  \ar@<-1ex>[l]
  \ar[r]
  \ar@<1ex>[l]
  &
C'_{2}
  \ar@<-1ex>[l]
  \ar[r]
  \ar@<1ex>[l]
  &
\cdots
  \ar@<-1ex>[l]
  \ar[r]
  \ar@<1ex>[l]
  &
C'_{n}
  \ar@<-1ex>[l]
  \ar@<1ex>[l]
}}
\]

As before, it follows that, for every \(0\leq k\leq m\), 
\[
F\bigr|^{C'_{m}}_{C_{m}}\bigr|^{(C'_{m})_{k}}_{(C_{m})_{k}}=F\bigr|^{C'_{k}}_{C_{k}}.
\]

This completes the proof.
\end{proof}

\begin{corollary}\label{CUPwCat}
The ordered pair \((\w\mathsf{Cat}, (\mathrm{Gd}^{(n)}\circ U^{(n,\omega)})_{n\in\omega})\) satisfies that, for every \(m,n\in\omega\) with \(m\leq n\),
\[
\mathrm{Gd}^{(m)}\circ U^{(m,\omega)}
=
U^{(m,n)}_{\mathrm{MS}}\circ\mathrm{Gd}^{(n)}\circ U^{(n,\omega)}.
\]
Therefore, we will denote by \(\mathrm{Gd}^{(\omega)}\) the functor \(\langle\mathrm{Gd}^{(n)}\circ U^{(n,\omega)}\rangle_{n\in\omega}\) from \(\w\mathsf{Cat}\) to \(\w\mathsf{Cat}^{\mathrm{MS}}\) defined in Proposition~\ref{PUnivPropMS}.
\end{corollary}

\begin{proof}
For every \(m,n\in\omega\) with \(m\leq n\), the following chain of equalities holds
\begin{align*}
U^{(m,n)}_{\mathrm{MS}}\circ\mathrm{Gd}^{(n)}\circ U^{(n,\omega)}
&=
\mathrm{Gd}^{(m)}\circ U^{(m,n)}\circ U^{(n,\omega)}
\tag{1}
\\
&=
\mathrm{Gd}^{(m)}\circ U^{(m,\omega)}.
\tag{2}
\end{align*}

The first equality follows from Proposition~\ref{PGradSq}; the second equality follows from Remark~\ref{RUmnProjSys}.
\end{proof}

\begin{remark}\label{RGdw}
Following the construction in the proof of Proposition~\ref{PUnivPropMS} and the definition of the graduation functors in Definition~\ref{DGrad} and \ref{DGradFun}, we now describe the explicit construction of the functor \(\mathrm{Gd}^{(\omega)}\).
\begin{enumerate}
\item
Given a single-sorted \(\omega\)-category \(\mathsf{C}=(C,(\#^{k},\mathrm{sc}^{k},\mathrm{tg}^{k})_{k\in\omega})\),
\[
\mathrm{Gd}^{(\omega)}(\mathsf{C})=
\left(
(C_{k})_{k\in\omega},
\left(
\circ^{(j)},\mathrm{sc}^{(j,k)},\mathrm{tg}^{(j,k)},\mathrm{i}^{(k,j)}
\right)_{(j,k)\in\coprod_{k\in\omega}k}
\right)
\]
is the many-sorted \(\omega\)-category where, for every \(k\in\omega\),
\[
C_{k}=\{f\in C\mid f=\mathrm{sc}^{k}(f)\},
\]
is the set of \(k\)-cells of \(\mathsf{C}\) and, for every \(0\leq j<k\), the structural operations are given by
\begin{align*}
\circ^{(j)}&=\#^{j}\big|^{C_{k}}_{C_{k}\times C_{k}};
&\mathrm{i}^{(k,j)}&=\mathrm{in}^{C_{j},C_{k}};
\\
\mathrm{sc}^{(j,k)}&=\mathrm{sc}^{j}\big|^{C_{j}}_{C_{k}};
&\mathrm{tg}^{(j,k)}&=\mathrm{tg}^{j}\big|^{C_{j}}_{C_{k}}.
\end{align*}
\item
Moreover, given two single-sorted \(\omega\)-categories \(\mathsf{C}\) and \(\mathsf{C}'\) and \(F:\mathsf{C}\mor\mathsf{C}'\) a single-sorted \(\omega\)-functor, \(\mathrm{Gd} ^{(\omega)}(F)=(F\big|^{C'_{k}}_{C_{k}})_{k\in\omega}\).
\end{enumerate}
\end{remark}

\section{
\texorpdfstring
{From many-sorted to single-sorted \(\omega\)-categories}
{From many-sorted to single-sorted Omega-categories}
}

In this section we associate to a many-sorted $\omega$-category a single-sorted $\omega$-category by means of the notion of \(\omega\)-unification. Moreover, we also associate to a many-sorted $\omega$-functor between many-sorted $\omega$-categories a single-sorted $\omega$-functor between the associated single-sorted $\omega$-categories. We end this section by proving that the above construction is functorial.

\begin{definition}\label{D[Cw]}
Let $\mathsf{C} = ((C_{k})_{k\in\omega},\zeta)$ be a many-sorted $\omega$-category  with 
$
\zeta = (\zeta_{j,k})_{(j,k)\in \coprod_{k\in\omega}k}
$ 
and, for every $(j,k)\in \coprod_{k\in\omega}k$,
$$
\zeta_{j,k} =  \left(
\circ^{(j)},\mathrm{sc}^{(j,k)},\mathrm{tg}^{(j,k)},\mathrm{i}^{(k,j)}
\right).
$$
Then, the set of \emph{\(\omega\)-cells} of \(\mathsf{C}\), denoted by \(C_{\omega}\) is,
\[\textstyle
C_{\omega}=\coprod_{k\in\omega}C_{k}.
\]
In particular, a pair \((f,j)\) is a \(\omega\)-cell, i.e., \((f,j)\in C_{\omega}\), if \(j\in\omega\) and \(f\in C_{j}\). We will say that \(f\) is the \emph{representative} and \(j\) is the \emph{level} of the \(\omega\)-cell \((f,j)\). For every \(j\in\omega\), we let \(\mathrm{in}^{C_{j}}\) stand for the inclusion mapping from \(C_{j}\) to \(C_{\omega}\) that assigns to every element \(f\) in \(C_{j}\) the element \((f,j)\) in \(C_{\omega}\).
\[
\mathrm{in}^{C_{j}}:C_{j}\mor C_{\omega}.
\]

Moreover, we let \(\Theta\) be the binary relation on \(C_{\omega}\) defined as \(((f,j),(g,k))\in\Theta\) if
\begin{enumerate}
\item
\(j<k\) and \(g=\mathrm{i}^{(k,j)}(f)\); or
\item
\(j=k\) and \(f=g\); or
\item
\(j>k\) and \(f=\mathrm{i}^{(j,k)}(g)\).
\end{enumerate}
It follows that \(\Theta\) is an equivalence relation on \(C_{\omega}\). Thus, we will write \([C_{\omega}]\) in place of the quotient set \(C_{\omega}/\Theta\). Moreover, for every \(\omega\)-cell \((f,j)\) in \(C_{\omega}\), we let \([f]_{j}\) stand for the equivalence class \([(f,j)]_{\Theta}\). The elements \([f]_{j}\) in \([C_{\omega}]\) will be called \emph{\(\omega\)-cell classes} and they will be denoted by \([f]_{j}\), \([g]_{k}\), \([h]_{l}\), ... . We let \(\mathrm{pr}^{\Theta}\) stand for the projection mapping from \(C_{\omega}\) to \([C_{\omega}]\) that assigns to a \(\omega\)-cell \((f,j)\) in \(C_{\omega}\) the \(\omega\)-cell class \([f]_{j}\).
\[
\mathrm{pr}^{\Theta}:C_{\omega}\mor[C_{\omega}].
\]
Moreover, for every \(j\in\omega\), we let \(\mathrm{in}^{[C_{j}]}\) stand for the mapping from \(C_{j}\) to \([C_{\omega}]\) given by the composition \(\mathrm{pr}^{\Theta}\circ\mathrm{in}^{C_{j}}\), that is, \(\mathrm{in}^{[C_{j}]}\) sends an element \(f\) in \(C_{j}\) to the \(\omega\)-cell class \([f]_{j}\) in \([C_{\omega}]\).
\[
\mathrm{in}^{[C_{j}]}:C_{j}\mor [C_{\omega}].
\]
Moreover, for every \(\omega\)-cell \((f,j)\) in \(C_{\omega}\), its \emph{support}, denoted \(\mathrm{supp}(f,j)\), is the set 
\[
\{i\in\omega\mid\exists x\in C_{i}, ((f,j),(x,i))\in\Theta\}.
\]
By definition of the equivalence relation \(\Theta\), for every \(i\in\mathrm{supp}(f,j)\), there exist a unique element \(x\) in \(C_{i}\) such that \(((f,j),(x,i))\in\Theta\). Moreover, \(\mathrm{supp}(f,j)\) is non-empty since \(j\) trivially belongs to \(\mathrm{supp}(f,j)\). For every \(\omega\)-cell \((f,j)\) in \(C_{\omega}\), \(\min\mathrm{supp}(f,j)\) will be called the \emph{nature} of \((f,j)\) in \(\mathsf{C}\) and we will denote it by \(\mathrm{nat}_{\mathsf{C}}(f,j)\). Moreover, its associated element of \(C_{\mathrm{nat}_{\mathsf{C}}(f,j)}\) the \emph{natural representative} of \((f,j)\). This definitions trivially extends to \(\omega\)-cell classes and we will also write \(\mathrm{supp}([f]_{j})\).

Finally, let \(\mathsf{C}=(C,(\xi_{k})_{k\in\omega})\) be a single-sorted \(\omega\)-category. Note that, for every \(f\in C\), by item (\(\omega\)C2) in Definition~\ref{DnCatSS}, the set
\[
\{i\in\omega\mid f=\mathrm{sc}^{i}(f)\},
\]
is non-empty. Therefore, its minimum will be called the \emph{nature} of \(f\) in \(\mathsf{C}\) and we will denote it by \(\mathrm{nat}_{\mathsf{C}}(f)\).
\end{definition}

\begin{remark}
The set \([C_{\omega}]\) is the inductive limit of the inductive system \((C_{(\bigcdot)},\mathrm{i}^{(\bigcdot^{+},\bigcdot)})\) where \(C_{(\bigcdot)}\) is the family of sets \((C_{k})_{k\in\omega}\) and \(\mathrm{i}^{(\bigcdot^{+},\bigcdot)}\) is the family of mappings \((\mathrm{i}^{(k+1,k)})_{k\in\omega}\). Thus, the following universal property follows.

Let \(L\) be a set and \((f_{i})_{i\in\omega}\) be a \(\omega\)-family of mappings such that, for every \(i\in\omega\), \(f_{i}\) is a mapping from \(C_{i}\) to \(L\). If, for every \(0\leq j<k\), \(f_{j}=f_{k}\circ\mathrm{i}^{(k,j)}\), then there exists a unique mapping \([f_{i}]_{i\in\omega}\) from \([C_{\omega}]\) to \(L\) such that, for every \(j\in\omega\), \(f_{j}=[f_{i}]_{i\in\omega}\circ\mathrm{in}^{[C_{j}]}\). This situation is depicted in Figure~\ref{FUnivProp[Cw]}.
\begin{figure}
\[
\xymatrix{
C_{j}
  \ar[rr]^{\mathrm{i}^{(k,j)}}
  \ar[rd]_{\mathrm{in}^{[C_{j}]}}
  \ar@/_2pc/[rdd]_{f_{j}}
&&
C_{k}
  \ar[ld]^{\mathrm{in}^{[C_{k}]}}
  \ar@/^2pc/[ldd]^{f_{k}}
\\
&
[C_{\omega}]
  \ar[d]^{[f_{i}]_{i\in\omega}}
\\
&
L
}
\]
\caption{Universal property of \([C_{\omega}]\)}
\label{FUnivProp[Cw]}
\end{figure}
\end{remark}

\begin{definition}\label{DScTg}
Let $\mathsf{C} = ((C_{k})_{k\in\omega},\zeta)$ be a many-sorted $\omega$-category  with 
$
\zeta = (\zeta_{j,k})_{(j,k)\in \coprod_{k\in\omega}k}
$ 
and, for every $(j,k)\in \coprod_{k\in\omega}k$,
$$
\zeta_{j,k} =  \left(
\circ^{(j)},\mathrm{sc}^{(j,k)},\mathrm{tg}^{(j,k)},\mathrm{i}^{(k,j)}
\right).
$$
Let \(i\in\omega\). We let \(\mathrm{sc}^{i}\) and \(\mathrm{tg}^{i}\) stand for the assignments from \([C_{\omega}]\) to \([C_{\omega}]\) defined, for every \(\omega\)-cell class \([f]_{j}\in [C_{\omega]}\) as
\begin{align*}
\mathrm{sc}^{i}\left(
[f]_{j}
\right)
&=
\begin{cases}
[f]_{j}&\mbox{if \(j\leq i\)};
\\
[\mathrm{sc}^{(i,j)}(f)]_{i}&\mbox{if \(i<f\)}.
\end{cases}
&
\mathrm{tg}^{i}\left(
[f]_{j}
\right)
&=
\begin{cases}
[f]_{j}&\mbox{if \(j\leq i\)};
\\
[\mathrm{tg}^{(i,j)}(f)]_{i}&\mbox{if \(i<j\)}.
\end{cases}
\end{align*}
\end{definition}

\begin{proposition}\label{PScTg}
Let $\mathsf{C} = ((C_{k})_{k\in\omega},\zeta)$ be a many-sorted $\omega$-category  with 
$
\zeta = (\zeta_{j,k})_{(j,k)\in \coprod_{k\in\omega}k}
$ 
and, for every $(j,k)\in \coprod_{k\in\omega}k$,
$$
\zeta_{j,k} =  \left(
\circ^{(j)},\mathrm{sc}^{(j,k)},\mathrm{tg}^{(j,k)},\mathrm{i}^{(k,j)}
\right).
$$
For every \(i\in\omega\), \(\mathrm{sc}^{i}\) and \(\mathrm{tg}^{i}\) are (well-defined) unary operations on \([C_{\omega}]\). That is, for every \([f]_{j}\) and \([g]_{k}\) \(\omega\)-cell classes in \([C_{\omega}]\), if \([f]_{j}=[g]_{k}\), then
\begin{align*}
\mathrm{sc}^{i}([f]_{j})&=\mathrm{sc}^{i}([g]_{k})
&
\mathrm{tg}^{i}([f]_{j})&=\mathrm{tg}^{i}([g]_{k}).
\end{align*}
\end{proposition}

\begin{proof}
Let \(i\in\omega\). Let us first note that it suffices to prove that, for every \(\omega\)-cell class \([f]_{j}\) in \([C_{\omega}]\), if \([f_{0}]_{j_{0}}\) is its natural representation, then
\[
\mathrm{sc}^{i}([f]_{j})
=
\mathrm{sc}^{i}\left([f_{0}]_{j_{0}}\right).
\tag{\(\star\)}
\]



Let \([f]_{j}\) be a \(\omega\)-cell class in \([C_{\omega}]\) and \([f_{0}]_{j_{0}}\) its natural representation. Therefore, according to the definition of \(\Theta\), either (A) \(j=j_{0}\) and \(f=f_{0}\) or (B) \(j_{0}<j\) and \(f=\mathrm{i}^{(j,j_{0})}(f_{0})\). Since the first case is trivial, let us consider the second one.

We want to prove \((\star)\). 

If (B), that is, if \(j_{0}<j\) and \(f=\mathrm{i}^{(j,j_{0})}(f_{0})\), then we distinguish the following cases (1) \(j_{0}<j\leq i\); (2) \(j_{0}< i<j\); (3) \(j_{0}=i<j\); and (4) \(i<j_{0}<j\).

If (1), i.e., if \(j_{0}< j\leq i\), then unraveling the definition of \(\mathrm{sc}^{i}\) presented in Definition~\ref{DScTg} it follows that
\[
\mathrm{sc}^{i}([f_{0}]_{j_{0}})=
[f_{0}]_{j_{0}}=
\left[\mathrm{i}^{(j,j_{0})}(f_{0})\right]_{j}=
\mathrm{sc}^{i}\left(\left[\mathrm{i}^{(j,j_{0})}(f_{0})\right]_{j}\right)
=
\mathrm{sc}^{i}\left(\left[f\right]_{j}\right)
.
\]

If (2), i.e., if \(j_{0}<i<j\) then the following chain of equalities holds
\allowdisplaybreaks
\begin{align*}
\mathrm{sc}^{i}\left(\left[f\right]_{j}\right)
&=
\mathrm{sc}^{i}\left(\left[\mathrm{i}^{(j,j_{0})}(f_{0})\right]_{j}\right)
\tag{1}
\\
&=
\left[\mathrm{sc}^{(i,j)}\left(\mathrm{i}^{(j,j_{0})}\left(f_{0}\right)\right)\right]_{i}
\tag{2}
\\
&=
\left[\mathrm{sc}^{(i,j)}\left(\mathrm{i}^{(j,i)}\left(\mathrm{i}^{(i,j_{0})}\left(f_{0}\right)\right)\right)\right]_{i}
\tag{3}
\\
&=
\left[\mathrm{i}^{(i,j_{0})}(f_{0})\right]_{i}
\tag{4}
\\
&=
[f_{0}]_{j_{0}}
\tag{5}
\\
&=
\mathrm{sc}^{i}([f_{0}]_{j_{0}}).
\tag{6}
\end{align*}

The first equality follows by replacing equals by equals; the second equality unravels the definition of \(\mathrm{sc}^{i}\) presented in Definition~\ref{DScTg} since \(i<j\); the third equality follows from item~(MS4) in Definition~\ref{DnCatMS}; the fourth equality follows from item~(MS5) in Definition~\ref{DnCatMS}; the fifth equality follows from the definition of \(\Theta\) presented in Definition~\ref{D[Cw]}; finally, the last equality recovers the definition of \(\mathrm{sc}^{i}\) since \(j_{0}<i\).

If (3), i.e., if \(j_{0}=i<j\), then the following chain of equalities holds
\allowdisplaybreaks
\begin{align*}
\mathrm{sc}^{i}\left(\left[f\right]_{j}\right)
&=
\mathrm{sc}^{j_{0}}\left(\left[\mathrm{i}^{(j,j_{0})}(f_{0})\right]_{j}\right)
\tag{1}
\\
&=
\left[\mathrm{sc}^{(j_{0},j)}\left(\mathrm{i}^{(j,j_{0})}\left(f_{0}\right)\right)\right]_{j_{0}}
\tag{2}
\\
&=
[f_{0}]_{j_{0}}
\tag{3}
\\
&=
\mathrm{sc}^{j_{0}}([f_{0}]_{j_{0}}).
\tag{4}
\end{align*}

The first equality follows by replacing equals by equals; the second equality unravels the definition of \(\mathrm{sc}^{j_{0}}\) presented in Definition~\ref{DScTg} since \(j_{0}<j\); the third equality follows from item~(MS5) in Definition~\ref{DnCatMS}; finally, the last equality recovers the definition of \(\mathrm{sc}^{j_{0}}\).

If (4), then the following chain of equalities holds
\begin{align*}
\mathrm{sc}^{i}\left(\left[f\right]_{j}\right)
&=
\mathrm{sc}^{i}\left(\left[\mathrm{i}^{(j,j_{0})}(f_{0})\right]_{j}\right)
\tag{1}
\\
&=
\left[\mathrm{sc}^{(i,j)}\left(\mathrm{i}^{(j,j_{0})}\left(f_{0}\right)\right)\right]_{i}
\tag{2}
\\
&=
\left[\mathrm{sc}^{(i,j_{0})}\left(\mathrm{sc}^{(j_{0},j)}\left(\mathrm{i}^{(j,j_{0})}\left(f_{0}\right)\right)\right)\right]_{i}
\tag{3}
\\
&=
\left[\mathrm{sc}^{(i,j_{0})}(f_{0})\right]_{i}
\tag{4}
\\
&=
\mathrm{sc}^{i}([f_{0}]_{j_{0}}).
\tag{5}
\end{align*}

The first equality follows by replacing equals by equals; the second equality unravels the definition of \(\mathrm{sc}^{i}\) presented in Definition~\ref{DScTg} since \(i<j\); the third equality follows from item~(MS1) in Definition~\ref{DnCatMS}; the fourth equality follows from item~(MS5) in Definition~\ref{DnCatMS}; finally, the last equality recovers the definition of \(\mathrm{sc}^{i}\) since \(i<j_{0}\).

By a similar argument follows that \(\mathrm{tg}^{i}\) is a unary operation on \([C_{\omega}]\).

This concludes the proof.
\end{proof}

\begin{definition}\label{DComp}
Let $\mathsf{C} = ((C_{k})_{k\in\omega},\zeta)$ be a many-sorted $\omega$-category  with 
$
\zeta = (\zeta_{j,k})_{(j,k)\in \coprod_{k\in\omega}k}
$ 
and, for every $(j,k)\in \coprod_{k\in\omega}k$,
$$
\zeta_{j,k} =  \left(
\circ^{(j)},\mathrm{sc}^{(j,k)},\mathrm{tg}^{(j,k)},\mathrm{i}^{(k,j)}
\right).
$$

Let \(i\in\omega\). We let \(\#^{i}\) stand for the partial binary operation from \([C_{\omega}]\times[C_{\omega}]\) to \([C_{\omega}]\) defined, for every \(\omega\)-cell classes \([f]_{j}\) and \([g]_{k}\) in \([C_{\omega}]\) with \(\mathrm{sc}^{i}([f]_{j})=\mathrm{tg}^{i}([g]_{k})\), as follows.
\begin{enumerate}
\item
If \(j\leq i\), then \(\mathrm{sc}^{i}([f]_{j})=[f]_{j}\) and \([f]_{j}\) will act as an \(i\)-cell identity. Therefore, we define \( [f]_{j}\#^{i}[g]_{k}=[g]_{k}\).
\item
If \(k\leq i\), then \(\mathrm{tg}^{i}([g]_{k})=[g]_{k}\) and \([g]_{k}\) will act as an \(i\)-cell identity. Therefore, we define \([f]_{j}\#^{i}[g]_{k}=[g]_{k}\).
\item
If \(i<j\) and \(i<k\), then we distinguish the following cases.
\begin{enumerate}
\item
If \(k<j\), then \([g]_{k}=[\mathrm{i}^{(j,k)}(g)]_{j}\). Let us note that the following chain of equalities holds
\begin{align*}
\left[\mathrm{sc}^{(i,j)}(f)\right]_{i}
&=
\mathrm{sc}^{i}([f]_{j})
\tag{1}
\\
&=
\mathrm{tg}^{i}([g]_{k})
\tag{2}
\\
&=
\mathrm{tg}^{i}\left(\left[\mathrm{i}^{(j,k)}(g)\right]_{j}\right)
\tag{3}
\\
&=
\left[\mathrm{tg}^{(i,j)}(\mathrm{i}^{(j,k)}(g))\right]_{i}.
\tag{4}
\end{align*}

The first equality follows from the definition of \(\mathrm{sc}^{i}\) presented in Definition~\ref{DScTg} since \(i<j\); the second equality is the initial assumption; the third equality follows from Proposition~\ref{PScTg}; finally, the last equality unravels the definition of \(\mathrm{tg}^{i}\) presented in Definition~\ref{DScTg} since \(i<j\).

Thus, \(\mathrm{sc}^{(i,j)}(f)=\mathrm{tg}^{(i,j)}(\mathrm{i}^{(j,k)}(g))\) and we define 
\[
[f]_{j}\#^{i}[g]_{k}=[f\circ^{(i)}\mathrm{i}^{(j,k)}(g)]_{j}.
\]
\item
Similarly, if \(j<k\), then \([f]_{j}=[\mathrm{i}^{(k,j)}(g)]_{k}\) and we define 
\[
[f]_{j}\#^{i}[g]_{k}=[\mathrm{i}^{(k,j)}(f)\circ^{(i)}g]_{k}.
\]
\item
Finally, if \(j=k\), then \(\mathrm{sc}^{i}([f]_{j})=\left[\mathrm{sc}^{(i,j)}(f)\right]_{i}\) and \(\mathrm{tg}^{i}([g]_{k})=\left[\mathrm{tg}^{(i,k)}(g)\right]_{i}\). Thus, \(\mathrm{sc}^{(i,j)}(f)=\mathrm{tg}^{(i,k)}(g)\) and we define 
\[
[f]_{j}\#^{i}[g]_{k}=[f\circ^{(i)}g]_{j}.
\]
\end{enumerate}
\end{enumerate}

Let us note that, if it is the case that both \(j\leq i\) and \(k\leq i\), then
\[
[f]_{j}=\mathrm{sc}^{i}([f]_{j})=\mathrm{tg}^{i}([g]_{k})=[g]_{k}
\]
and the definition is sound.
\end{definition}

\begin{proposition}\label{PComp}
Let $\mathsf{C} = ((C_{k})_{k\in\omega},\zeta)$ be a many-sorted $\omega$-category  with 
$
\zeta = (\zeta_{j,k})_{(j,k)\in \coprod_{k\in\omega}k}
$ 
and, for every $(j,k)\in \coprod_{k\in\omega}k$,
$$
\zeta_{j,k} =  \left(
\circ^{(j)},\mathrm{sc}^{(j,k)},\mathrm{tg}^{(j,k)},\mathrm{i}^{(k,j)}
\right).
$$
For every \(i\in\omega\), \(\#^{i}\) is a (well-defined) binary partial operation on \([C_{\omega}]\). That is, for every \([f]_{j},[f']_{j'},[g]_{k}\) and \([g']_{k'}\) \(\omega\)-cell classes in \([C_{\omega}]\) with \([f]_{j}=[f']_{j'}\) and \([g]_{k}=[g']_{k'}\) and \(\mathrm{sc}^{i}([f]_{j})=\mathrm{tg}^{i}([g]_{k})\) (which by Proposition~\ref{PScTg}, implies that \(\mathrm{sc}^{i}([f']_{j'})=\mathrm{tg}^{i}([g']_{k'})\)), then
\[
[f]_{j}\#^{i}[g]_{k}=[f']_{j'}\#^{i}[g']_{k'}.
\]
\end{proposition}

\begin{proof}
Let \(i\in\omega\). Note that it suffices to show that, for every \(\omega\)-cell classes \([f]_{j}\) and \([g]_{k}\) in \([C_{\omega}]\), if \([f_{0}]_{j_{0}}\) and \([g_{0}]_{k_{0}}\) are their natural representatives, respectively, then
\begin{align*}
[f]_{j}\#^{i}[g]_{k}&=[f_{0}]_{j_{0}}\#^{i}[g]_{k}
&\mbox{and}&&
[f]_{j}\#^{i}[g]_{k}&=[f]_{j}\#^{i}[g_{0}]_{k_{0}}
\tag{\(\star\)}
\end{align*}

We will only prove to prove the right-hand side of \((\star)\).  The left-hand side equation of \((\star)\) is done similarly. 

Let \([f]_{j}\) and \([g]_{k}\) be \(\omega\)-cell classes in \([C_{\omega}]\) and let \([g_{0}]_{k_{0}}\) be the natural representation of \([g]_{k}\). 

Since \([g]_{k}=[g_{0}]_{k_{0}}\), according to the definition of \(\Theta\), either (A) \(k=k_{0}\) and \(g=g_{0}\) or (B) \(k_{0}<k\) and \(g=\mathrm{i}^{(k,k_{0})}(g_{0})\). Since the first case is trivial, let us consider the second one.

We will consider a series of cases. The structure of the proof is as follows. First we unravel the cases which define the \(\#^{i}\)-composite \([f]_{j}\#^{i}[g]_{k}\). After that, in the respective cases, we unravel the subcases which define the \(\#^{i}\)-composite \([f]_{j}\#^{i}[g_{0}]_{k_{0}}\) to prove the result.

Before proceeding with the proof, we prove the following claim about many-sorted \(\omega\)-categories.

\begin{claim}\label{CId}
Let \(k\leq i<k'\), \(g\in C_{k}\) and \(f\in C_{k'}\). If \(\mathrm{sc}^{(i,k')}(f)=\mathrm{tg}^{(i,k')}(\mathrm{i}^{(k',k)}(g))\), then
\[
f\circ^{(i)}\mathrm{i}^{(k',k)}(g)=f.
\]
\end{claim}

The following chain of equalities holds
\begin{align*}
\mathrm{sc}^{(i,k')}(f)
&=
\mathrm{tg}^{(i,k')}\left(\mathrm{i}^{(k',k)}\left(g\right)\right)
\tag{1}
\\
&=
\mathrm{tg}^{(i,k')}\left(\mathrm{i}^{(k',i)}\left(\mathrm{i}^{(i,k)}\left(g\right)\right)\right)
\tag{2}
\\
&=
\mathrm{i}^{(i,k)}\left(g\right).
\tag{3}
\end{align*}

The first equality follows from assumption; the second equality follows from item~(MS4) in Definition~\ref{DnCatMS}; finally, the last equality follows from item~(MS5) in Definition~\ref{DnCatMS}.

Thus, the following chain of equalities holds
\begin{align*}
f\circ^{(i)}\mathrm{i}^{(k',k)}(g)
&=
f\circ^{(i)}\mathrm{i}^{(k',i)}\left(\mathrm{i}^{(i,k)}\left(g\right)\right)
\tag{1}
\\
&=
f\circ^{(i)}\mathrm{i}^{(k',i)}\left(\mathrm{sc}^{(i,k')}\left(f\right)\right)
\tag{2}
\\
&=
f.
\tag{3}
\end{align*}

The first equality follows from item~(MS4) in Definition~\ref{DnCatMS}; the second equality follows from the previously stated chain of equalities; finally, the last equality follows from item~(MS6) in Definition~\ref{DnCatMS}.

This completes the proof of Claim~\ref{CId}.

We are now in a position to prove \((\star)\). Following Definition~\ref{DComp} to determine the \(\#^{i}\)-composite \([f]_{j}\#^{i}[g]_{k}\), we consider the following cases (1) \(j\leq i\); (2) \(k\leq i\); or (3) \(i<j\) and \(i<k\).

\textsf{Case (1).}
If \(j\leq i\), by definition of \(\#^{i}\),
\[
[f]_{j}\#^{i}[g]_{k}=[g]_{k}=[g_{0}]_{k_{0}}=[f]_{j}\#^{i}[g_{0}]_{k_{0}}.
\]
Thus, \((\star)\) holds in Case 1.

\textsf{Case (2).}
If \(k\leq i\), since \(k_{0}<k\), then \(k_{0}< i\). Thus, by definition of \(\#^{i}\),
\[
[f]_{j}\#^{i}[g]_{k}=[f]_{j}=[f]_{j}\#^{i}[g_{0}]_{k_{0}}.
\]
Thus, \((\star)\) holds in Case 2.

\textsf{Case (3).}
If both \(i<j\) and \(i<k\), then we distinguish the following cases (3A) \(j<k\); (3B) \(k<j\); or (3C) \(j=k\).

\textsf{Case (3A).}
If \(i<j\), \(i<k\) and \(j<k\), then
\[\label{EComp1}
[f]_{j}\#^{i}[g]_{k}
=
\left[\mathrm{i}^{(k,j)}(f)\circ^{(i)}g\right]_{k}
=
\left[\mathrm{i}^{(k,j)}(f)\circ^{(i)}\mathrm{i}^{(k,k_{0})}(g_{0})\right]_{k}.
\tag{Comp1}
\]

Now, following Definition~\ref{DComp} to determine the \(\#^{i}\)-composite \([f]_{j}\#^{i}[g_{0}]_{k_{0}}\), we consider the following cases (3A1) \(j\leq i\); (3A2) \(k_{0}\leq i\); or (3A3) \(i<j\) and \(i<k_{0}\).

\textsf{Case (3A1).}
If \(i<j\), \(i<k\), \(j<k\) and \(j\leq i\), then we immediately get to a contradiction since \(i<j\leq i\). Thus, Case 3A1 is impossible.

\textsf{Case (3A2).}
If \(i<j\), \(i<k\), \(j<k\) and \(k_{0}\leq i\), then the following chain of equalities holds
\begin{align*}
[f]_{j}\#^{i}[g_{0}]_{k_{0}} 
&=
[f]_{j}
\tag{1}
\\
&=
\left[\mathrm{i}^{(k,j)}(f)\right]_{k}
\tag{2}
\\
&=
\left[\mathrm{i}^{(k,j)}(f)\circ\mathrm{i}^{(k,k_{0})}(g_{0})\right]_{k}
\tag{3}
\\
&=
[f]_{j}\#^{i}[g]_{k}.
\tag{4}
\end{align*}

The first equality unravels the definition of \(\#^{i}\) presented in Definition~\ref{DComp} since \(k\leq i\); the second equality follows from the definition of the relation \(\Theta\) presented in Definition~\ref{D[Cw]}; the third equality follows from Claim~\ref{CId}; finally, the last equality follows from Equation~(\ref{EComp1}).

Thus, \((\star)\) holds in Case 3A2.

\textsf{Case (3A3).}
If \(i<j\), \(i<k\), \(j<k\) and \(i<k_{0}\), then we distinguish the following cases (3A3A) \(j<k_{0}\); or (3A3B) \(k_{0}<j\); or (3A3C) \(j=k_{0}\).

\textsf{Case (3A3A).}
If \(i<j\), \(i<k\), \(j<k\), \(i<k_{0}\) and \(j<k_{0}\), then the following chain of equalities holds
\begin{align*}
[f]_{j}\#^{i}[g_{0}]_{k_{0}}
&=
\left[\mathrm{i}^{(k_{0},j)}(f)\circ^{(i)}g_{0}\right]_{k_{0}}
\tag{1}
\\
&=
\left[\mathrm{i}^{(k,k_{0})}\left(\mathrm{i}^{(k_{0},j)}(f)\circ^{(i)}g_{0}\right)\right]_{k}
\tag{2}
\\
&=
\left[\mathrm{i}^{(k,k_{0})}\left(\mathrm{i}^{(k_{0},j)}(f)\right)\circ^{(i)}\mathrm{i}^{(k,k_{0})}(g_{0})\right]_{k}
\tag{3}
\\
&=
\left[\mathrm{i}^{(k,j)}(f)\circ^{(i)}\mathrm{i}^{(k,k_{0})}(g_{0})\right]_{k}
\tag{4}
\\
&=
[f]_{j}\#^{i}[g]_{k}.
\tag{5}
\end{align*}

The first equality unravels the definition of \(\#^{i}\) presented in Definition~\ref{DComp} since \(i<j<k_{0}\); the second equality follows from the definition of \(\Theta\) presented in Definition~\ref{D[Cw]}; the third equality follows from item~(MS8) in Definition~\ref{DnCatMS} since \(i<k_{0}<k\); the fourth equality follows from item~(MS4) in Definition~\ref{DnCatMS}; finally, the last equality follows from Equation~(\ref{EComp1}).

Thus, \((\star)\) holds in Case 3A3A.

\textsf{Case (3A3B).}
If \(i<j\), \(i<k\), \(j<k\), \(i<k_{0}\) and \(k_{0}<j\), then the following chain of equalities holds
\begin{align*}
[f]_{j}\#^{i}[g_{0}]_{k_{0}}
&=
\left[f\circ^{(i)}\mathrm{i}^{(j,k_{0})}(g_{0})\right]_{j}
\tag{1}
\\
&=
\left[\mathrm{i}^{(k,j)}\left(f\circ^{(i)}\mathrm{i}^{(j,k_{0})}(g_{0})\right)\right]_{k}
\tag{2}
\\
&=
\left[\mathrm{i}^{(k,j)}\left(f\right)\circ^{(i)}\mathrm{i}^{(k,j)}\left(\mathrm{i}^{(j,k_{0})}(g_{0})\right)\right]_{k}
\tag{3}
\\
&=
\left[\mathrm{i}^{(k,j)}\left(f\right)\circ^{(i)}\mathrm{i}^{(k,k_{0})}(g_{0})\right]_{k}
\tag{4}
\\
&=
[f]_{j}\#^{i}[g]_{k}.
\tag{5}
\end{align*}

The first equality unravels the definition of \(\#^{i}\) presented in Definition~\ref{DComp} since \(i<k_{0}<j\); the second equality follows from the definition of \(\Theta\) presented in Definition~\ref{D[Cw]}; the third equality follows from item~(MS8) in Definition~\ref{DnCatMS} since \(i<j<k\); the fourth equality follows from item~(MS4) in Definition~\ref{DnCatMS}; finally, the last equality follows from Equation~(\ref{EComp1}).

Thus, \((\star)\) holds in Case 3A3B.

\textsf{Case (3A3C).}
If \(i<j\), \(i<k\), \(j<k\), \(i<k_{0}\) and \(j=k_{0}\), then the following chain of equalities holds
\begin{align*}
[f]_{j}\#^{i}[g_{0}]_{k_{0}}
&=
\left[f\circ^{(i)}g_{0}\right]_{j}
\tag{1}
\\
&=
\left[\mathrm{i}^{(k,j)}\left(f\circ^{(i)}g_{0}\right)\right]_{k}
\tag{2}
\\
&=
\left[\mathrm{i}^{(k,J)}\left(f\right)\circ^{(i)}\mathrm{i}^{(k,j)}\left(g_{0}\right)\right]_{k}
\tag{3}
\\
&=
\left[\mathrm{i}^{(k,j)}\left(f\right)\circ^{(i)}\mathrm{i}^{(k,k_{0})}(g_{0})\right]_{k}
\tag{4}
\\
&=
[f]_{j}\#^{i}[g]_{k}.
\tag{5}
\end{align*}

The first equality unravels the definition of \(\#^{i}\) presented in Definition~\ref{DComp} since \(j=k_{0}\); the second equality follows from the definition of \(\Theta\) presented in Definition~\ref{D[Cw]}; the third equality follows from item~(MS8) in Definition~\ref{DnCatMS} since \(i<j<k\); the fourth equality follows since \(j=k_{0}\) in Definition~\ref{DnCatMS}; finally, the last equality follows from Equation~(\ref{EComp1}).

Thus, \((\star)\) holds in case 3A3C.

This completes the proof for Case 3A.

\textsf{Case (3B).}
If \(i<j\), \(i<k\) and \(k<j\), then
\begin{align*}
[f]_{j}\#^{i}[g]_{k}
&=
\left[f\circ^{(i)}\mathrm{i}^{(j,k)}\left(\mathrm{i}^{(k,k_{0})}(g_{0})\right)\right]_{j}
\\\label{EComp2}
&=
\left[f\circ^{(i)}\mathrm{i}^{(j,k_{0})}(g_{0})\right]_{j}.
\tag{Comp2}
\end{align*}
The first equality unravels the definition of \(\#^{i}\); the second equality follows from item~(MS4) in Definition~\ref{DnCatMS}.

Now, following Definition~\ref{DComp} to determine the \(\#^{i}\)-composite \([f]_{j}\#^{i}[g_{0}]_{k_{0}}\), we consider the following cases (3B1) \(j\leq i\); (3B2) \(k_{0}\leq i\); or (3B3) \(i<j\) and \(i<k_{0}\). 

\textsf{Case (3B1).}
If \(i<j\), \(i<k\), \(k<j\) and \(j\leq i\), then we immediately get to a contradiction since \(i<j\leq i\). Thus, Case (3B1) is impossible.

\textsf{Case (3B2).}
If \(i<j\), \(i<k\), \(k<j\) and \(k_{0}\leq i\), then the following chain of equalities holds
\begin{align*}
[f]_{j}\#^{i}[g_{0}]_{k_{0}} 
&=
[f]_{j}
\tag{1}
\\
&=
\left[f\circ\mathrm{i}^{(j,k_{0})}(g_{0})\right]_{j}
\tag{2}
\\
&=
[f]_{j}\#^{i}[g]_{k}.
\tag{3}
\end{align*}

The first equality unravels the definition of \(\#^{i}\) presented in Definition~\ref{DComp} since \(k_{0}\leq i\); the second equality follows from Claim~\ref{CId}; finally, the last equality follows from Equation~(\ref{EComp2}).

Thus, \((\star)\) holds in Case 3B2.

\textsf{Case (3B3).}
If \(i<j\), \(i<k\), \(k<j\) and \(i<k_{0}\), since \(k<j\) and \(k_{0}<k\), it follows that \(k_{0}<j\). Therefore, the following chain of equalities holds
\[
[f]_{j}\#^{i}[g_{0}]_{k_{0}}
=
\left[f\circ^{(i)}\mathrm{i}^{(j,k_{0})}(g_{0})\right]_{j}
=
[f]_{j}\#^{i}[g]_{k}.
\]

Thus, \((\star)\) holds in Case 3B3.

This completes the proof for Case 3B.

\textsf{Case (3C).}
If \(i<j\), \(i<k\) and \(j=k\), then
\[\label{EComp3}
[f]_{j}\#^{i}[g]_{k}
=
\left[f\circ^{(i)}\mathrm{i}^{(k,k_{0})}(g_{0})\right]_{k}.
\tag{Comp3}
\]

Now, following Definition~\ref{DComp} to determine the \(\#^{i}\)-composite \([f]_{j}\#^{i}[g_{0}]_{k_{0}}\), we consider the following cases (3C1) \(j\leq i\); (3C2) \(k_{0}\leq i\); or (3C3) \(i\leq j\) and \(i<k_{0}\). 

\textsf{Case (3C1).}
If \(i<j\), \(i<k\), \(k=j\) and \(j\leq i\), then we immediately get to a contradiction since \(i<j\leq i\). Thus, Case 3C1 is impossible.

\textsf{Case (3C2).}
If \(i<j\), \(i<k\), \(k=j\) and \(k_{0}\leq i\), then the following chain of equalities holds
\begin{align*}
[f]_{j}\#^{i}[g_{0}]_{k_{0}} 
&=
[f]_{j}
\tag{1}
\\
&=
\left[f\circ\mathrm{i}^{(k,k_{0})}(g_{0})\right]_{k}
\tag{2}
\\
&=
[f]_{j}\#^{i}[g]_{k}.
\tag{3}
\end{align*}

The first equality unravels the definition of \(\#^{i}\) presented in Definition~\ref{DComp} since \(k_{0}\leq i\); the second equality follows from Claim~\ref{CId} since \(j=k\); finally, the last equality follows from Equation~(\ref{EComp3}).

Thus, \((\star)\) holds in Case 3C2.

\textsf{Case (3C3).}
If \(i<j\), \(i<k\), \(k=j\) and \(i<k_{0}\), since \(j=k\) and \(k_{0}<k\), it follows that \(k_{0}<j\). Therefore, the following chain of equalities holds
\[
[f]_{j}\#^{i}[g_{0}]_{k_{0}}
=
\left[f\circ^{(i)}\mathrm{i}^{(j,k_{0})}(g_{0})\right]_{j}
=
[f]_{j}\#^{i}[g]_{k}.
\]

Thus, \((\star)\) holds in Case 3C3.

This completes the proof for Case 3C.

This completes the proof for Case 3.

Thus, let \([f]_{j},[f']_{j'},[g]_{k},[g']_{k'}\in[C_{\omega}]\) with \([f]_{j}=[f']_{j'}\), \([g]_{k}=[g']_{k'}\) and \(\mathrm{sc}^{i}([f]_{j})=\mathrm{tg}^{i}([g]_{k})\). Let \([f_{0}]_{j_{0}}\) be the natural representation of \([f]_{j}=[f']_{j'}\) and \([g_{0}]_{k_{0}}\) be the natural representation of \([g]_{k}=[g']_{k'}\). Then, the following chain of equalities holds
\begin{align*}
[f]_{j}\#^{i}[g]_{k}
&=
[f_{0}]_{j_{0}}\#^{i}[g]_{k}
\tag{1}
\\
&=
[f']_{j'}\#^{i}[g]_{k}
\tag{2}
\\
&=
[f']_{j'}\#^{i}[g_{0}]_{k_{0}}
\tag{3}
\\
&=
[f']_{j'}\#^{i}[g']_{k'}
\tag{4}
\end{align*}

This completes the proof.
\end{proof}

\begin{definition}\label{DUfwCat}
Let $\mathsf{C} = ((C_{k})_{k\in\omega},\zeta)$ be a many-sorted $\omega$-category  with 
$
\zeta = (\zeta_{j,k})_{(j,k)\in \coprod_{k\in\omega}k}
$ 
and, for every $(j,k)\in \coprod_{k\in\omega}k$,
$$
\zeta_{j,k} =  \left(
\circ^{(j)},\mathrm{sc}^{(j,k)},\mathrm{tg}^{(j,k)},\mathrm{i}^{(k,j)}
\right).
$$
We will call the ordered pair
\[
\left(
[C_{\omega}],
\left(\#^{i},\mathrm{sc}^{i},\mathrm{tg}^{i}\right)_{i\in\omega}
\right)
\]
the \emph{unification of \(\mathsf{C}\)} and we will denote it by \(\mathrm{Uf}^{(\omega)}(\mathsf{C})\).
\end{definition}

We are now in a position to prove that, for a many-sorted \(\omega\)-category \(\mathsf{C}\), its unification \(\mathrm{Uf}^{(\omega)}(\mathsf{C})\) is a single-sorted \(\omega\)-category.

\begin{proposition}\label{PUfwCat}
Let $\mathsf{C} = ((C_{k})_{k\in\omega},\zeta)$ be a many-sorted $\omega$-category  with 
$
\zeta = (\zeta_{j,k})_{(j,k)\in \coprod_{k\in\omega}k}
$ 
and, for every $(j,k)\in \coprod_{k\in\omega}k$,
$$
\zeta_{j,k} =  \left(
\circ^{(j)},\mathrm{sc}^{(j,k)},\mathrm{tg}^{(j,k)},\mathrm{i}^{(k,j)}
\right).
$$
Then \(\mathrm{Uf}^{(\omega)}(\mathsf{C})\) is a single-sorted \(\omega\)-category.
\end{proposition}

\begin{proof}
We want to prove that \((\xi_{i})_{i\in\omega}=(\#^{i},\mathrm{sc}^{i},\mathrm{tg}^{i})_{i\in\omega}\) as described in Definitions~\ref{D[Cw]} and \ref{DComp} is a structure of single-sorted \(\omega\)-category on \([C_{\omega}]\). Thus, we must check that it satisfies the conditions stated in Definition~\ref{DnCatSS}.

We begin by proving that, for \(k\in n\), \(([C_{\omega}],(\#^{i},\mathrm{sc}^{i},\mathrm{tg}^{i}))\) is a single-sorted category. To this end, we consider the different items stated in Definition~\ref{DCat}.

\textsf{(C1)}
Let \([f]_{j}\in[C_{\omega}]\). Consider \(m=\max\{i,j\}+1\). The following chain of equalities holds
\begin{align*}
\mathrm{sc}^{i}\left(\mathrm{sc}^{i}\left(\left[f\right]_{j}\right)\right)
&=
\mathrm{sc}^{i}\left(\mathrm{sc}^{i}\left(\left[\mathrm{i}^{(m,j)}(f)\right]_{m}\right)\right)
\tag{1}
\\
&=
\mathrm{sc}^{i}\left(\left[\mathrm{sc}^{(i,m)}\left(\mathrm{i}^{(m,j)}(f)\right)\right]_{i}\right)
\tag{2}
\\
&=
\left[\mathrm{sc}^{(i,m)}\left(\mathrm{i}^{(m,j)}(f)\right)\right]_{i}
\tag{3}
\\
&=
\mathrm{sc}^{i}\left(\left[\mathrm{i}^{(m,j)}(f)\right]_{m}\right)
\tag{4}
\\
&=
\mathrm{sc}^{i}\left(\left[f\right]_{j}\right).
\tag{5}
\end{align*}

The first equality follows from Proposition~\ref{PScTg}; the second and third equality unravels the definition of \(\mathrm{sc}^{i}\) presented in Definition~\ref{DScTg} since \(i<m\) and \(i\leq i\); the fourth equality recovers the definition of \(\mathrm{sc}^{i}\); finally, the last equality follows again from Proposition~\ref{PScTg}. 

By a similar argument we also have that
\begin{align*}
\mathrm{sc}^{i}\circ\mathrm{tg}^{i}=\mathrm{tg}^{i};&
&\mathrm{tg}^{i}\circ\mathrm{sc}^{i}=\mathrm{sc}^{i};&
&\mathrm{tg}^{i}\circ\mathrm{tg}^{i}=\mathrm{tg}^{i}.
\end{align*}

\textsf{(C2)}
For every \([f]_{j},[g]_{k}\in[C_{\omega}]\), the equivalence
\begin{align*}
[f]_{j}\#^{i}[g]_{k}&
&\mbox{if and only if}&&
\mathrm{sc}^{i}([f]_{j})&=\mathrm{tg}^{i}([g]_{k})
\end{align*}
follows from the definition of \(\#^{i}\) presented in Definition~\ref{DComp}.

\textsf{(C3)}
Let \([f]_{j},[g]_{k}\in[C_{\omega}]\) with \(\mathrm{sc}^{i}\left([f]_{j}\right)=\mathrm{tg}^{i}\left([g]_{k}\right)\). Consider \(m=\max\{i,j,k\}+1\). The following chain of equalities holds
\allowdisplaybreaks
\begin{align*}
\mathrm{sc}^{i}\left(\left[f\right]_{j}\#^{i}\left[g\right]_{k}\right)
&=
\mathrm{sc}^{i}\left(\left[\mathrm{i}^{(m,j)}(f)\right]_{m}\#^{i}\left[\mathrm{i}^{(m,k)}(g)\right]_{m}\right)
\tag{1}
\\
&=
\mathrm{sc}^{i}\left(\left[\mathrm{i}^{(m,j)}(f)\circ^{(i)}\mathrm{i}^{(m,k)}(g)\right]_{m}\right)
\tag{2}
\\
&=
\left[\mathrm{sc}^{(i,m)}\left(\mathrm{i}^{(m,j)}(f)\circ^{(i)}\mathrm{i}^{(m,k)}(g)\right)\right]_{i}
\tag{3}
\\
&=
\left[\mathrm{sc}^{(i,m)}\left(\mathrm{i}^{(m,k)}(g)\right)\right]_{i}
\tag{4}
\\
&=
\mathrm{sc}^{i}\left(\left[\mathrm{i}^{(m,k)}(g)\right]_{m}\right)
\tag{5}
\\
&=
\mathrm{sc}^{i}\left(\left[g\right]_{k}\right).
\tag{6}
\end{align*}

The first equality follows from Proposition~\ref{PComp}; the second equality unravels the definition of \(\#^{i}\) presented in Definition~\ref{DComp} since \(i<m\); the third equality unravels the definition of \(\mathrm{sc}^{i}\) presented in Definition~\ref{DScTg} since \(i<m\); The fourth equality follows from item~(MS2) in Definition~\ref{DnCatMS}; the fifth equality recovers the definition of \(\mathrm{sc}^{i}\); finally, the last equality follows from Proposition~\ref{PScTg}.

By a similar argument we also have that
\[
\mathrm{tg}^{i}\left([f]_{j}\#^{i}[g]_{k}\right)=\mathrm{tg}^{i}\left([f]_{j}\right).
\]

\textsf{(C4)}
Let \([f]_{j}\in[C_{\omega}]\). Consider \(m=\max\{i,j\}+1\). The following chain of equalities holds
\begin{align*}
[f]_{j}\#^{i}\mathrm{sc}^{i}\left(\left[f\right]_{j}\right)
&=
[f]_{j}\#^{i}\mathrm{sc}^{i}\left(\left[\mathrm{i}^{(m,j)}(f)\right]_{m}\right)
\tag{1}
\\
&=
[f]_{j}\#^{i}\left[\mathrm{sc}^{(i,m)}\left(\mathrm{i}^{(m,j)}(f)\right)\right]_{i}
\tag{2}
\\
&=
\left[\mathrm{i}^{(m,j)}(f)\right]_{m}\#^{i}\left[\mathrm{i}^{(m,i)}\left(\mathrm{sc}^{(i,m)}\left(\mathrm{i}^{(m,j)}(f)\right)\right)\right]_{m}
\tag{3}
\\
&=
\left[\mathrm{i}^{(m,j)}(f)\circ^{(i)}\mathrm{i}^{(m,i)}\left(\mathrm{sc}^{(i,m)}\left(\mathrm{i}^{(m,j)}(f)\right)\right)\right]_{m}
\tag{4}
\\
&=
\left[\mathrm{i}^{(m,j)}(f)\right]_{m}
\tag{5}
\\
&=
[f]_{j}.
\tag{6}
\end{align*}

The first equality follows from Proposition~\ref{PScTg}; the second equality unravels the definition of \(\mathrm{sc}^{i}\) presented in Definition~\ref{DScTg} since \(i<m\); the third equality follows from Proposition~\ref{PComp}; the fourth equality unravels the definition of \(\#^{i}\) presented in Definition~\ref{DComp}; the fifth equality follows from item~(MS8) in Definition~\ref{DnCatMS}; finally, the last equality follows from the definition of the equivalence relation \(\Theta\).

By a similar argument we also have that
\[
\mathrm{tg}^{i}\left(\left[f\right]_{j}\right)\#^{i}\left[f\right]_{j}
=
\left[f\right]_{j}.
\]

\textsf{(C5)}
Let \([f]_{j},[g]_{k},[h]_{l}\in[C_{\omega}]\) with \(\mathrm{sc}^{i}([f]_{j})=\mathrm{tg}^{i}([g]_{k})\) and \(\mathrm{sc}^{i}([g]_{k})=\mathrm{tg}^{i}([h]_{l})\). Consider \(m=\max\{i,j,k,l\}+1\). The following chain of equalities holds.
\begin{align*}
\left[f\right]_{j}\#^{i}\left(\left[g\right]_{k}\#^{i}\left[h\right]_{l}\right)
&=
\left[\mathrm{i}^{(m,j)}(f)\right]_{m}\#^{i}\left(\left[\mathrm{i}^{(m,k)}(g)\right]_{m}\#^{i}\left[\mathrm{i}^{(m,l)}(h)\right]_{m}\right)
\tag{1}
\\
&=
\left[\mathrm{i}^{(m,j)}(f)\circ^{(i)}\left(\mathrm{i}^{(m,k)}(g)\circ^{(i)}\mathrm{i}^{(m,l)}(h)\right)\right]_{m}
\tag{2}
\\
&=
\left[\left(\mathrm{i}^{(m,j)}(f)\circ^{(i)}\mathrm{i}^{(m,k)}(g)\right)\circ^{(i)}\mathrm{i}^{(m,l)}(h)\right]_{m}
\tag{3}
\\
&=
\left(\left[\mathrm{i}^{(m,j)}(f)\right]_{m}\#^{i}\left[\mathrm{i}^{(m,k)}(g)\right]_{m}\right)\#^{i}\left[\mathrm{i}^{(m,l)}(h)\right]_{m}
\tag{4}
\\
&=
\left(\left[f\right]_{j}\#^{i}\left[g\right]_{k}\right)\#^{i}\left[h\right]_{l}.
\tag{5}
\end{align*}

The first equality follows from Proposition~\ref{PComp}; the second equality unravels the definition of \(\#^{i}\) presented in Definition~\ref{DComp} since \(i<m\); the third equality follows from item~(MS7) in Definition~\ref{DnCatMS}; the fourth equality recovers the definition of \(\#^{i}\); finally, the last equality follows again from Proposition~\ref{PComp}.

It follows that, for every \(i\in\omega\), \(([C_{\omega}], (\#^{i},\mathrm{sc}^{i},\mathrm{tg}^{i}))\) is a single-sorted category.

According to Definition~\ref{DnCatSS}, now we have to prove that, for every \(n,m\in\omega\) with \(n<m\), the ordered pair \(([C_{\omega}], (\#^{i},\mathrm{sc}^{i},\mathrm{tg}^{i})_{i\in\{n,m\}})\) is a \(2\)-category. To this end, we consider the different items stated in Definition~\ref{D2Cat}.

\textsf{(2C1)}
Let \([f]_{j}\in[C_{\omega}]\). Consider \(p=\max\{m,j\}+1\). Then the following chain of equalities holds
\begin{align*}
\mathrm{sc}^{n}\left(\mathrm{sc}^{m}\left([f]_{j}\right)\right)
&=
\mathrm{sc}^{n}\left(\mathrm{sc}^{m}\left(\left[\mathrm{i}^{(p,j)}(f)\right]_{p}\right)\right)
\tag{1}
\\
&=
\mathrm{sc}^{n}\left(\left[\mathrm{sc}^{(m,p)}\left(\mathrm{i}^{(p,j)}(f)\right)\right]_{m}\right)
\tag{2}
\\
&=
\left[\mathrm{sc}^{(n,m)}\left(\mathrm{sc}^{(m,p)}\left(\mathrm{i}^{(p,j)}(f)\right)\right)\right]_{n}
\tag{3}
\\
&=
\left[\mathrm{sc}^{(n,p)}\left(\mathrm{i}^{(p,j)}(f)\right)\right]_{n}
\tag{4}
\\
&=
\mathrm{sc}^{n}\left(\left[\mathrm{i}^{(p,j)}(f)\right]_{p}\right)
\tag{5}
\\
&=
\mathrm{sc}^{n}\left([f]_{j}\right).
\tag{6}
\end{align*}

The first equality follows from Proposition~\ref{PScTg}; the second and third equalities unravel the definition of \(\mathrm{sc}^{m}\) and \(\mathrm{sc}^{n}\), respectively, presented in Definition~\ref{DScTg} since \(n<m<p\); the fourth equality follows from item~(MS1) in Definition~\ref{DnCatMS}; the fifth equality recovers the definition of \(\mathrm{sc}^{n}\); finally, the last equality follows from Proposition~\ref{PScTg}.

By a similar argument we also have that 
\[
\mathrm{sc}^{n}(\mathrm{tg}^{m}([f]_{j}))=\mathrm{sc}^{n}([f]_{j}).
\]

On the other hand, the following chain of equalities holds
\begin{align*}
\mathrm{sc}^{m}\left(\mathrm{sc}^{n}\left([f]_{j}\right)\right)
&=
\mathrm{sc}^{m}\left(\mathrm{sc}^{n}\left(\left[\mathrm{i}^{(p,j)}(f)\right]_{p}\right)\right)
\tag{1}
\\
&=
\mathrm{sc}^{m}\left(\left[\mathrm{sc}^{(n,p)}\left(\mathrm{i}^{(p,j)}(f)\right)\right]_{n}\right)
\tag{2}
\\
&=
\left[\mathrm{sc}^{(n,p)}\left(\mathrm{i}^{(p,j)}(f)\right)\right]_{n}
\tag{3}
\\
&=
\mathrm{sc}^{n}\left(\left[\mathrm{i}^{(p,j)}(f)\right]_{p}\right)
\tag{4}
\\
&=
\mathrm{sc}^{n}\left([f]_{j}\right).
\tag{5}
\end{align*}

The first equality follows from Proposition~\ref{PScTg}; the second and third equalities unravel the definitions of \(\mathrm{sc}^{n}\) and \(\mathrm{sc}^{m}\), respectively, presented in Definition~\ref{DScTg} since \(n<p\) and \(n< m\); the fourth equality recovers the definition of \(\mathrm{sc}^{n}\); finally, the last equality follows from Proposition~\ref{PScTg}.

By a similar argument we also have that
\[
\mathrm{tg}^{m}\circ \mathrm{tg}^{n}
=
\mathrm{tg}^{n}
=
\mathrm{tg}^{n}\circ \mathrm{tg}^{m}
=
\mathrm{tg}^{n}\circ \mathrm{sc}^{m}.
\]

\textsf{(2C2)}
Let \([f]_{j},[g]_{k}\in[C_{\omega}]\) with \(\mathrm{sc}^{n}([f]_{j})=\mathrm{tg}^{n}([g_{k}])\). Consider \(p=\max\{j,k,m\}+1\). Then the following chain of equalities holds
\begin{align*}
\mathrm{sc}^{m}\left([f]_{j}\#^{n}[g]_{k}\right)
&=
\mathrm{sc}^{m}\left(\left[i^{(p,j)}(f)\right]_{p}\#^{n}\left[i^{(p,k)}(g)\right]_{p}\right)
\tag{1}
\\
&=
\mathrm{sc}^{m}\left(\left[\mathrm{i}^{(p,j)}(f)\circ^{(n)}\mathrm{i}^{(p,k)}(g)\right]_{p}\right)
\tag{2}
\\
&=
\left[\mathrm{sc}^{(m,p)}\left(\mathrm{i}^{(p,j)}(f)\circ^{(n)}\mathrm{i}^{(p,k)}(g)\right)\right]_{m}
\tag{3}
\\
&=
\left[\mathrm{sc}^{(m,p)}\left(\mathrm{i}^{(p,j)}(f)\right)\circ^{(n)}\mathrm{sc}^{(m,p)}\left(\mathrm{i}^{(p,k)}(g)\right)\right]_{m}
\tag{4}
\\
&=
\left[\mathrm{sc}^{(m,p)}\left(\mathrm{i}^{(p,j)}(f)\right)\right]_{m}\#^{n}\left[\mathrm{sc}^{(m,p)}\left(\mathrm{i}^{(p,k)}(g)\right)\right]_{m}
\tag{5}
\\
&=
\mathrm{sc}^{m}\left(\left[\mathrm{i}^{(p,j)}(f)\right]_{p}\right)\#^{n}\mathrm{sc}^{m}\left(\left[\mathrm{i}^{(p,k)}(g)\right]_{p}\right)
\tag{6}
\\
&=
\mathrm{sc}^{m}(\left[f\right]_{j})\#^{n}\mathrm{sc}^{m}\left(\left[g\right]_{k}\right).
\tag{7}
\end{align*}

The first equality follows from Proposition~\ref{PComp}; the second and third equalities unravel the definitions of \(\#^{n}\) and \(\mathrm{sc}^{m}\), respectively, presented in Definitions~\ref{DComp} and \ref{DScTg} since \(n<m<p\); the fourth equality follows from item~(MS3) of Definition~\ref{DnCatMS}; the fifth and sixth equalities recover the definitions of \(\#^{n}\) and \(\mathrm{sc}^{m}\), respectively; finally, the last equality follows from Proposition~\ref{PScTg}.

\textsf{(2C3)}
Let \([f]_{j},[f']_{j'},[g]_{k},[g']_{k'}\in[C_{\omega}]\) with
\begin{align*}
\mathrm{sc}^{m}\left([f]_{j}\right)&=\mathrm{tg}^{m}\left([f']_{j'}\right),
&\mathrm{sc}^{m}\left([g]_{k}\right)&=\mathrm{tg}^{m}\left([g']_{k'}\right),
\\
\mathrm{sc}^{n}\left([f]_{j}\right)&=\mathrm{tg}^{n}\left([g]_{k}\right),
&\mathrm{sc}^{n}\left([f']_{j'}\right)&=\mathrm{tg}^{n}\left([g']_{k'}\right).
\end{align*}
Consider \(p=\max\{j,j',k,k',m\}+1\). Then, the following chain of equalities holds
\begin{align*}
&\left([f]_{j}\#^{n}[g]_{k}\right)
\#^{m}
\left([f']_{j'}\#^{n}[g']_{k'}\right)
\\
&=
\left(\left[i^{(p,j)}(f)\right]_{p}\#^{n}\left[i^{(p,k)}(g)\right]_{p}\right)
\#^{m}
\left(\left[i^{(p,j')}(f')\right]_{p}\#^{n}\left[i^{(p,k')}(g')\right]_{p}\right)
\tag{1}
\\
&=
\left(\left[i^{(p,j)}(f)\circ^{(n)} i^{(p,k)}(g)\right]_{p}\right)
\#^{m}
\left(\left[i^{(p,j')}(f')\circ^{(n)} i^{(p,k')}(g')\right]_{p}\right)
\tag{2}
\\
&=
\left[\left(i^{(p,j)}(f)\circ^{(n)} i^{(p,k)}(g)\right)
\circ^{(m)}
\left(i^{(p,j')}(f')\circ^{(n)} i^{(p,k')}(g')\right)\right]_{p}
\tag{3}
\\
&=
\left[\left(i^{(p,j)}(f)\circ^{(m)} i^{(p,j')}(f')\right)
\circ^{(n)}
\left(i^{(p,k)}(g)\circ^{(m)} i^{(p,k')}(g')\right)\right]_{p}
\tag{4}
\\
&=
\left(\left[i^{(p,j)}(f)\circ^{(m)} i^{(p,j')}(f')\right]_{p}\right)
\#^{n}
\left(\left[i^{(p,k)}(g)\circ^{(m)} i^{(p,k')}(g')\right]_{p}\right)
\tag{5}
\\
&=
\left(\left[i^{(p,j)}(f)\right]_{p}\#^{m}\left[i^{(p,j')}(f')\right]_{p}\right)
\#^{n}
\left(\left[i^{(p,k)}(g)\right]_{p}\#^{m}\left[i^{(p,k')}(g')\right]_{p}\right)
\tag{6}
\\
&=
\left([f]_{j}\#^{m}[f']_{j'}\right)
\#^{n}
\left([g]_{k}\#^{m}[g']_{k'}\right).
\tag{7}
\end{align*}

The first equality follows from Proposition~\ref{PComp}; the second and third equalities unravel the definitions of \(\#^{n}\) and \(\#^{m}\) presented in Definition~\ref{DComp}, respectively, presented in Definition~\ref{DComp}, the fourth equality follows from item~(MS9) of Definition~\ref{DnCatMS}; the fifth and sixth equalities recover the definitions of \(\#^{n}\) and \(\#^{m}\), respectively; finally, the last equality follows from Proposition~\ref{PComp}.

It follows that, for every \(n,m\in\omega\) with \(n<m\), \(([C_{\omega}], (\#^{i},\mathrm{sc}^{i},\mathrm{tg}^{i})_{i\in\{n,m\}})\) is a single-sorted 2-category, thus satisfying item \(\omega\)C1 in Definition~\ref{DnCatSS}.

Finally, Let us note that, for every \([f]_{j}\in[C_{\omega}]\), \([f]_{j}=\mathrm{sc}^{j}([f]_{j})\), thus satisfying item \(\omega\)C2 in Definition~\ref{DnCatSS}.

All in all, we can affirm that \(\mathrm{Uf}^{(\omega)}(\mathsf{C})\) is a single-sorted \(\omega\)-category.
\end{proof}

In the following definition we introduce the notion of unification of a many-sorted \(\omega\)-functor.

\begin{definition}\label{DUfwFun}
Let \(\mathsf{C}=((C_{k})_{k\in\omega},\zeta)\) with \(\zeta=(\zeta_{j,k})_{(j,k)\in\coprod_{k\in\omega}k}\) and, for every \((j,k)\in\coprod_{k\in\omega}k\),
\[
\zeta_{j,k}=\left(\circ^{(j)},\mathrm{sc}^{(j,k)},\mathrm{tg}^{(j,k)},\mathrm{i}^{(k,j)}\right)
\mbox{, and}
\]
\(\mathsf{C}'=((C'_{k})_{k\in\omega},\zeta')\) with \(\zeta'=(\zeta'_{j,k})_{(j,k)\in\coprod_{k\in\omega}k}\) and, for every \((j,k)\in\coprod_{k\in\omega}k\),
\[
\zeta'_{j,k}=\left(\circ'^{(j)},\mathrm{sc}'^{(j,k)},\mathrm{tg}'^{(j,k)},\mathrm{i}'^{(k,j)}\right)
\]
be two many-sorted \(\omega\)-categories, \(F:\mathsf{C}\mor\mathsf{C}'\) a many-sorted \(\omega\)-functor from \(\mathsf{C}\) to \(\mathsf{C}'\), and 
\begin{align*}
\mathrm{Uf}^{(\omega)}(\mathsf{C})&=\left([C_{\omega}], (\#^{i},\mathrm{sc}^{i},\mathrm{tg}^{i})_{i\in\omega}\right)
&
\mathrm{Uf}^{(\omega)}(\mathsf{C}')&=\left([C'_{\omega}], (\#'^{i},\mathrm{sc}'^{i},\mathrm{tg}'^{i})_{i\in\omega}\right)
\end{align*}
the unifications corresponding to \(\mathsf{C}\) and \(\mathsf{C}'\), respectively.

We will call the mapping from \([C_{\omega}]\) to \([C'_{\omega}]\) that sends every \([f]_{j}\in[C_{\omega}]\) to \([F_{j}(f)]_{j}\), the \emph{unification} of \(F\) and we will denote it by \(\mathrm{Uf}^{(\omega)}(F)\).
\end{definition}

\begin{proposition}\label{PUfwFunwd}
Let \(\mathsf{C}\) and \(\mathsf{C}'\) two many-sorted \(\omega\)-categories and \(F:\mathsf{C}\mor\mathsf{C}'\) a many-sorted \(\omega\)-functor from \(\mathsf{C}\) to \(\mathsf{C}'\). Then \(\mathrm{Uf}^{(\omega)}(F)\) is a well-defined mapping from \([C_{\omega}]\) to \([C'_{\omega}]\). That is, if \([f]_{j}\) and \([g]_{k}\) are \(\omega\)-cell classes in \([C_{\omega}]\), if \([f]_{j}=[g]_{k}\), then
\[
\mathrm{Uf}^{(\omega)}(F)([f]_{j})
=
\mathrm{Uf}^{(\omega)}(F)([g]_{k}).
\]
\end{proposition}

\begin{proof}
Note that it suffices to prove that, for every \(\omega\)-cell class \([f]_{j}\) in \([C_{\omega}]\), if \([f_{0}]_{j_{0}}\) is its natural representation, then
\[
\mathrm{Uf}^{(\omega)}(F)\left(\left[f\right]_{j}\right)
=
\mathrm{Uf}^{(\omega)}(F)\left(\left[f_{0}\right]_{j_{0}}\right).
\tag{\(\star\)}
\]

Let \([f]_{j}\) be a \(\omega\)-cell class in \([C_{\omega}]\) and \([f_{0}]_{j_{0}}\) its natural representation. Therefore, according to the definition of \(\Theta\), either (A) \(j=j_{0}\) and \(f=f_{0}\) or (B) \(j_{0}<j\) and \(f=\mathrm{i}^{(j,j_{0})}(f_{0})\). Since the first case is trivial let us consider the second one.

We want to prove (\(\star\)).

If (B), that is, if \(j_{0}<j\) and \(f=\mathrm{i}^{(j,j_{0})}(f_{0})\), then the following chain of equalities holds
\begin{align*}
\mathrm{Uf}^{(\omega)}(F)\left(\left[f\right]_{j}\right)
&=
\mathrm{Uf}^{(\omega)}(F)\left(\left[\mathrm{i}^{(j,j_{0})}(f_{0})\right]_{j}\right)
\tag{1}
\\
&=
\left[F_{j}(\mathrm{i}^{(j,j_{0})}(f_{0}))\right]_{j}
\tag{2}
\\
&=
\left[\mathrm{i}'^{(j,j_{0})}(F_{j_{0}}(f_{0}))\right]_{j}
\tag{3}
\\
&=
\left[F_{j_{0}}(f_{0})\right]_{j_{0}}
\tag{4}
\\
&=
\mathrm{Uf}^{(\omega)}(F)\left(\left[f_{0}\right]_{j_{0}}\right).
\tag{5}
\end{align*}

The first equality follows by replacing equals by equals; the second equality unravels the definition of \(\mathrm{Uf}^{(\omega)}(F)\) presented in Definition~\ref{DUfwFun}; the third equality follows from item~(MSii) in Definition~\ref{DnCatMS}; the fourth equality follows from the definition of \(\Theta\) presented in Definition~\ref{D[Cw]}; finally, the last equality recovers the definition of \(\mathrm{Uf}^{(\omega)}(F)\).
\end{proof}

In the following proposition we prove that, for every many-sorted \(\omega\)-functor \(F\), its unification \(\mathrm{Uf}^{(\omega)}(F)\) is a single-sorted \(\omega\)-functor.

\begin{proposition}\label{PUfwFun}
Let \(\mathsf{C}\) and \(\mathsf{C}'\) be two many-sorted \(\omega\)-categories and \(F:\mathsf{C}\mor\mathsf{C}'\) a many-sorted \(\omega\)-functor from \(\mathsf{C}\) to \(\mathsf{C}'\). Then \(\mathrm{Uf}^{(\omega)}(F)\) is a single-sorted \(\omega\)-functor from \(\mathsf{Uf}^{(\omega)}(\mathsf{C})\) to \(\mathrm{Uf}^{(\omega)}(\mathsf{C}')\).
\end{proposition}

\begin{proof}
It suffices to prove that the single-sorted mapping \(\mathrm{Uf}^{(\omega)}(F)\) described in Definition~\ref{DUfwFun} satisfies the conditions of Definition~\ref{DnCatSS}, thus defining a single-sorted \(\omega\)-functor from \(\mathrm{Uf}^{(\omega)}(\mathsf{C})\) to \(\mathrm{Uf}^{(\omega)}(\mathsf{C}')\).

Let \(i\in\omega\).

\textsf{(Ci)}
Let \([f]_{j}\in[C_{\omega}]\) and consider \(m=\max\{i,j\}+1\). The following chain of equalities holds
\begin{align*}
\mathrm{Uf}^{(\omega)}(F)\left(\mathrm{sc}^{i}\left(\left[f
\right]_{j}\right)\right)
&=
\mathrm{Uf}^{(\omega)}(F)\left(\mathrm{sc}^{i}\left(\left[\mathrm{i}^{(m,j)}(f
)\right]_{m}\right)\right)
\tag{1}
\\
&=
\mathrm{Uf}^{(\omega)}(F)\left(\left[\mathrm{sc}^{(i,m)}\left(\mathrm{i}^{(m,j)}(f)\right)\right]_{i}\right)
\tag{2}
\\
&=
\left[F_{i}\left(\mathrm{sc}^{(i,m)}\left(\mathrm{i}^{(m,j)}(f)\right)\right)\right]_{i}
\tag{3}
\\
&=
\left[\mathrm{sc}'^{(i,m)}\left(F_{m}\left(\mathrm{i}^{(m,j)}(f)\right)\right)\right]_{i}
\tag{4}
\\
&=
\mathrm{sc}'^{i}\left(\left[F_{m}\left(\mathrm{i}^{(m,j)}(f)\right)\right]_{m}\right)
\tag{5}
\\
&=
\mathrm{sc}'^{i}\left(\mathrm{Uf}^{(\omega)}(F)\left(\left[\mathrm{i}^{(m,j)}(f)\right]_{m}\right)\right)
\tag{6}
\\
&=
\mathrm{sc}'^{i}\left(\mathrm{Uf}^{(\omega)}(F)\left(\left[f\right]_{j}\right)\right).
\tag{7}
\end{align*}

The first equality follows from Proposition~\ref{PScTg}; the second and third equalities unravel the definitions of \(\mathrm{sc}^{i}\) and \(\mathrm{Uf}^{(\omega)}(F)\) presented in Definition~\ref{D[Cw]} and \ref{DUfwFun}, respectively; the fourth equality follows from item~(MSi) in Definition~\ref{DnCatMS}; the fifth and sixth equalities recover the definitions of \(\mathrm{sc}'^{i}\) and \(\mathrm{Uf}^{(\omega)}(F)\), respectively; finally, the last equality follows from Proposition~\ref{PUfwFunwd}.

By a similar argument we also have that
\[
\mathrm{Uf}^{(\omega)}(F)\left(\mathrm{tg}^{i}\left([f]_{j}\right)\right)
=
\mathrm{tg}'^{i}\left(\mathrm{Uf}^{(\omega)}(F)\left([f]_{j}\right)\right).
\]

\textsf{(Cii)}
Let \([f]_{j},[g]_{k}\in[C_{\omega}]\) with \(\mathrm{sc}^{i}\left([f]_{j}\right)=\mathrm{tg}^{i}\left([g]_{k}\right)\) and consider \(m=\max\{i,j,k\}+1\). The following chain of equalities holds
\allowdisplaybreaks
\begin{align*}
\mathrm{Uf}^{(\omega)}(F)\left([f]_{j}\#^{i}[g]_{k}\right)
&=
\mathrm{Uf}^{(\omega)}(F)\left(\left[\mathrm{i}^{(m,j)}(f)\right]_{m}\#^{i}\left[\mathrm{i}^{(m,k)}(g)\right]_{m}\right)
\tag{1}
\\
&=
\mathrm{Uf}^{(\omega)}(F)\left(\left[\mathrm{i}^{(m,j)}(f)\circ^{(i)}\mathrm{i}^{(m,k)}(g)\right]_{m}\right)
\tag{2}
\\
&=
\left[F_{m}\left(\mathrm{i}^{(m,j)}(f)\circ^{(i)}\mathrm{i}^{(m,k)}(g)\right)\right]_{m}
\tag{3}
\\
&=
\left[F_{m}\left(\mathrm{i}^{(m,j)}(f)\right)\circ'^{(i)}F_{m}\left(\mathrm{i}^{(m,k)}(g)\right)\right]_{m}
\tag{4}
\\
&=
\left[F_{m}\left(\mathrm{i}^{(m,j)}(f)\right)\right]\#'^{i}\left[F_{m}\left(\mathrm{i}^{(m,k)}(g)\right)\right]_{m}
\tag{5}
\\
&=
\mathrm{Uf}^{(\omega)}(F)\left(\left[\mathrm{i}^{(m,j)}(f)\right]_{m}\right) \#'^{i} \mathrm{Uf}^{(\omega)}(F)\left(\left[\mathrm{i}^{(m,k)}(g)\right]_{m}\right)
\tag{6}
\\
&=
\mathrm{Uf}^{(\omega)}(F)\left(\left[f\right]_{j}\right) \#'^{i} \mathrm{Uf}^{(\omega)}(F)\left(\left[g\right]_{k}\right).
\tag{7}
\end{align*}

The first equality follows from Proposition~\ref{PComp}; the second and third equalities unravel the definitions of \(\#^{i}\) and \(\mathrm{Uf}^{(\omega)}(F)\) presented in Definitions~\ref{DComp} and \ref{DUfwFun}, respectively; the fourth equality follows from item~(MSii) in Definition~\ref{DnCatMS}; the fifth and sixth equalities recover the definitions of \(\#'^{i}\) and \(\mathrm{Uf}^{(\omega)}(F)\), respectively; finally, the last equality follows from Proposition~\ref{PUfwFunwd}.

Therefore, \(\mathrm{Uf}^{(\omega)}(F)\) is a single-sorted \(\omega\)-functor from \(\mathrm{Uf}^{(\omega)}(\mathsf{C})\) to \(\mathrm{Uf}^{(\omega)}(\mathsf{C}')\).
\end{proof}

We are now in position to prove that the unification construction is functorial.

\begin{proposition}\label{PUfwFunctor}
\(\mathrm{Uf}^{(\omega)}:\w\mathsf{Cat}^{\mathrm{MS}}\mor\w\mathsf{Cat}\) is a functor.
\end{proposition}

\begin{proof}
Taking into account Propositions~\ref{PUfwCat} and \ref{PUfwFunctor} it suffices to prove that \(\mathrm{Uf}^{(\omega)}\) preserves many-sorted identity \(\omega\)-functors and compositions of many-sorted \(\omega\)-functors.

\textsf{Preservation of identities.}
Let \(\mathsf{C}\) be a many-sorted \(\omega\)-category and consider the many-sorted \(\omega\)-functor \(\mathrm{Id}^{\mathsf{C}}:\mathsf{C}\mor\mathsf{C}\), where \(\mathrm{Id}^{\mathsf{C}}=(\mathrm{Id}^{\mathrm{C}_{k}})_{k\in\omega}\). Then, for every \([f]_{j}\in[C_{\omega}]\), the following chain of equalities holds
\begin{align*}
\mathrm{Uf}^{(\omega)}(\mathrm{Id}^{\mathsf{C}})\left([f]_{j}\right)
&=
[\mathrm{Id}^{C_{j}}(f)]_{j}
\tag{1}
\\
&=
[f]_{j}
\tag{2}
\\
&=
\mathrm{Id}^{\mathrm{Uf}^{(\omega)}(\mathsf{C})}\left(\left[f\right]_{j}\right).
\tag{3}
\end{align*}

The first equality unravels Definition~\ref{DUfwFun} on the unification of \(\mathrm{Id}^{\mathsf{C}}\); finally, the second and third equalities follows by definition of the suitable identity mappings.

Thus, \(\mathrm{Uf}^{(\omega)}\) preserves identities.

\textsf{Preservation of compositions.}
Let \(F:\mathsf{D}\mor \mathsf{E}\) and \(G:\mathsf{C}\mor\mathsf{D}\) be two many-sorted \(\omega\)-functors between many-sorted \(\omega\)-categories, where \(G=(G_{k})_{k\in\omega}\) and \(F=(F_{k})_{k\in\omega}\). Then we want to check that
\[
\mathrm{Uf}^{(\omega)}(F\circ G)
=
\mathrm{Uf}^{(\omega)}(F)\circ\mathrm{Uf}^{(\omega)}(G).
\]

For every \([f]_{j}\in[C_{\omega}]\), the following chain of equalities holds
\begin{align*}
\mathrm{Uf}^{(\omega)}(F\circ G)\left(\left[f\right]_{j}\right)
&=
\mathrm{Uf}^{(\omega)}((F_{k}\circ G_{k})_{k\in\omega})\left(\left[f\right]_{j}\right)
\tag{1}
\\
&=
\left[F_{j}\left(G_{j}\left(f\right)\right)\right]_{j}
\tag{2}
\\
&=
\mathrm{Uf}^{(\omega)}(F)\left(\left[G_{j}\left(f\right)\right]_{j}\right)
\tag{3}
\\
&=
\mathrm{Uf}^{(\omega)}(F)\circ \mathrm{Uf}^{(\omega)}(G)\left(\left[f\right]_{j}\right).
\tag{4}
\end{align*}

The first equality unravels the definition of composition of many-sorted \(\omega\)-functors from Definition~\ref{DnCatMS}; the second equality unravels the definition of \(\omega\)-unification of a many-sorted \(\omega\)-functor presented in Definition~\ref{DUfwFun}; The third and fourth equalities recover the unification of the functors \(F\) and \(G\), respectively.

Thus, \(\mathrm{Uf}^{(\omega)}\) preserves compositions.

This completes the proof.
\end{proof}

\section{
\texorpdfstring
{On the relation between \(\omega\)-graduation and \(\omega\)-unification}
{On the relation between Omega-graduation and Omega-unification}
}

In this section, we begin by investigating the relation between the functors $\mathrm{Gd}^{(\omega)}$ and $\mathrm{Uf}^{(\omega)}$.
\begin{center}
\begin{tikzpicture}
[ACliment/.style={-{To [angle'=45, length=5.75pt, width=4pt, round]}}
, scale=1, 
AClimentD/.style={dogble eqgal sign distance,
-implies
}
]

\node[] (0R) at (0.5,0) [color=white] {};
\node[] (1L) at (2.3,0)  [color=white] {};

\node[] (C0) at (0,0) [] {$\w\mathsf{Cat}$};
\node[] (C1) at (3,0) [] {$\w\mathsf{Cat}^{\mathrm{MS}}$};

\draw[ACliment] (0R.north east) to node [above] {$\mathrm{Gd}^{(\omega)}$} (1L.north west);
\draw[ACliment] (1L.south west) to node [below] {$\mathrm{Uf}^{(\omega)}$} (0R.south east);
\end{tikzpicture}
\end{center}

We start by proving that the process of \(\omega\)-graduation followed by that of \(\omega\)-unification, is naturally isomorphic to the identity functor on the category of many-sorted \(\omega\)-categories.

\begin{proposition}\label{PUfwGdw} 
$\mathrm{Uf}^{(\omega)}\circ\mathrm{Gd}^{(\omega)}\cong\mathrm{Id}^{\w\mathsf{Cat}}$.
\end{proposition}

\begin{proof}
Let \(\mathsf{C}=(C,(\xi_{i})_{i\in\omega})\) with \((\xi_{i})_{i\in\omega}=(\#^{i},\mathrm{sc}^{i},\mathrm{tg}^{i})_{i\in\omega}\) be a single-sorted \(\omega\)-category and 
\[
\mathrm{Gd}^{(\omega)}(\mathsf{C})=
\left(
(C_{k})_{k\in\omega},
\left(
\circ^{(j)},\mathrm{sc}^{(j,k)},\mathrm{tg}^{(j,k)},\mathrm{i}^{(k,j)}
\right)_{(j,k)\in\coprod_{k\in\omega}k}
\right)
\]
the graduation of \(\mathsf{C}\). According to Remark~\ref{RGdw}, for every \(k\in\omega\),
\[
C_{k}=\{f\in C\mid f=\mathrm{sc}^{k}(f)\}.
\]

Moreover, for every \(0\leq j<k\), the structural operations of the graduation are given by
\allowdisplaybreaks
\begin{align*}
\circ^{(j)}&=\#^{j}\big|^{C_{k}}_{C_{k}\times C_{k}};
&\mathrm{i}^{(k,j)}&=\mathrm{in}^{C_{j},C_{k}};
\\
\mathrm{sc}^{(j,k)}&=\mathrm{sc}^{j}\big|^{C_{j}}_{C_{k}};
&\mathrm{tg}^{(j,k)}&=\mathrm{tg}^{j}\big|^{C_{j}}_{C_{k}}.
\end{align*}

Now, let us consider the single-sorted \(\omega\)-category given by the unification of the graduation of \(\mathsf{C}\), i.e.,
\[
\mathrm{Uf}^{(\omega)}\left(\mathrm{Gd}^{(\omega)}\left(\mathsf{C}\right)\right)=
\left(
[C_{\omega}],\left(\#'^{i},\mathrm{sc}'^{i},\mathrm{tg}'^{i}\right)_{i\in\omega}
\right),
\]
where \([C_{\omega}]\) is the quotient set \(C_{\omega}/\Theta\) presented in Definition~\ref{D[Cw]} and its structural operations are defined by cases in Definition~\ref{DScTg} and \ref{DComp}.

Recall that, according to Definition~\ref{D[Cw]}, for every \(f\in C\), its nature in \(\mathsf{C}\) is the minimum of the set
\[
\{i\in\omega \mid f=\mathrm{sc}^{i}(f)\},
\]
and it is denoted by \(\mathrm{nat}_{\mathsf{C}}(f)\). We let \(\varepsilon^{\mathsf{C}}\) stand for the mapping from \(C\) to \(C_{\omega}\) that assigns to every element \(f\in C\) the \(\omega\)-cell \((f,\mathrm{nat}(f))\). Moreover, we let \(\varepsilon^{[\mathsf{C}]}\) stand for the mapping from \(C\) to \([C_{\omega}]\) given by the composition \(\mathrm{pr}^{\Theta}\circ\varepsilon^{\mathsf{C}}\), i.e., \(\varepsilon^{[\mathsf{C}]}\) sends an element \(f\in C\) to the \(\omega\)-cell class \([f]_{\mathrm{nat}(f)}\).
\[
\varepsilon^{[\mathsf{C}]}:C\mor[C_{\omega}].
\]

Moreover, we let \(\eta^{\mathsf{C}}\) stand for the mapping from \(C_{\omega}\) to \(C\) that assigns to every \(\omega\)-cell \((f,j)\) in \(C_{\omega}\) the element \(f \in C_{j}\subseteq C\).

\begin{claim}\label{CUfwGdwA}
For every \(\omega\)-cells \((f,j)\) and \((g,k)\) in \(C_{\omega}\), if \(\left((f,j),(g,k)\right)\in\Theta\), then \(\eta^{\mathsf{C}}(f,j)=\eta^{\mathsf{C}}(g,k)\).
\end{claim}

We want to prove that \(\eta^{\mathsf{C}}(f,j)=\eta^{\mathsf{C}}(g,k)\), that is, \(f=g\). We distinguish the following cases (1) \(j<k\); or (2) \(j=k\); or (3) \(k<j\).

If (1), i.e., if \(j<k\), then
\[
g=
\mathrm{i}^{(k,j)}(f)=
\mathrm{in}^{C_{k},C_{j}}(f)=
f.
\]

If (2), i.e., if \(j=k\), then, by the definition of \(\Theta\), it follows that \(f=g\).

If (3), i.e., if \(k<j\), then
\[
f=
\mathrm{i}^{(j,k)}(g)=
\mathrm{in}^{C_{j},C_{k}}(g)=
g.
\]

Claim~\ref{CUfwGdwA} follows.

Before proceeding any further, let us remark that Claim~\ref{CUfwGdwA} allows us to give a concrete description of any equivalence class in \(\mathrm{Uf}^{(\omega)}(\mathrm{Gd}^{(\omega)}(\mathsf{C}))\). In particular, if \([f]_{j}\in [C_{\omega}]\), then
\[
[f]_{j}=\{(f,i)\in C_{\omega}\mid \mathrm{nat}_{\mathsf{C}}(f)\leq i\}.
\]
Therefore, any two representations of the equivalence class \([f]_{j}\) share the same representative.

We let \(\eta^{[\mathsf{C}]}=\left(\eta^{\mathsf{C}}\right)^{\natural}\) stand for the mapping obtained by the universal property of the quotient as the unique mapping from \([C_{\omega}]\) to \(C\) such that \(\left(\eta^{\mathsf{C}}\right)^{\natural}\circ\mathrm{pr}^{\Theta}=\eta^{\mathsf{C}}\). Thus, \(\eta^{[\mathsf{C}]}\) is the mapping from \([C_{\omega}]\) to \(C\) that assigns to a \(\omega\)-cell class \([f]_{j}\) in \([C_{\omega}]\) the element \(\eta^{\mathsf{C}}(f,j)\), i.e., the element \(f\in C\). That \(\eta^{[\mathsf{C}]}\) is well-defined follows from Claim~\ref{CUfwGdwA}.
\[
\eta^{[\mathsf{C}]}:[C_{\omega}]\mor C.
\]

Our next aim is to prove that, for every single-sorted \(\omega\)-category \(\mathsf{C}\), \(\varepsilon^{[\mathsf{C}]}\) is a single-sorted \(\omega\)-functor.

\begin{claim}\label{CUfwGdwB}
\(\varepsilon^{[\mathsf{C}]}\) is a single-sorted \(\omega\)-functor from \(\mathsf{C}\) to \(\mathrm{Uf}^{(\omega)}(\mathrm{Gd}^{(\omega)}(\mathsf{C}))\).
\end{claim}

It suffices to prove that \(\varepsilon^{[\mathsf{C}]}\) satisfies the conditions of Definition~\ref{DnCatSS}.

Let \(i\in\omega\). We want to prove that \(\varepsilon^{[\mathsf{C}]}\) is a single-sorted functor from \((C, (\#^{i},\mathrm{sc}^{i},\mathrm{tg}^{i}))\) to \(([C_{\omega}], (\#'^{i},\mathrm{sc}'^{i},\mathrm{tg}'^{i}))\). To this end, we consider the different items stated in Definition~\ref{DCat}.

\textsf{(Ci)}
Let \(f\) be an element of \(C\). Consider \(m=\max\{\mathrm{nat}_{\mathsf{C}}(f), i\}+1\). The following chain of equalities holds
\allowdisplaybreaks
\begin{align*}
\mathrm{sc}'^{i}\left(\varepsilon^{[\mathsf{C}]}(f)\right)
&=
\mathrm{sc}'^{i}\left(\left[f\right]_{\mathrm{nat}_{\mathsf{C}}(f)}\right)
\tag{1}
\\
&=
\mathrm{sc}'^{i}\left(\left[f\right]_{m}\right)
\tag{2}
\\
&=
\left[\mathrm{sc}^{(i,m)}(f)\right]_{i}
\tag{3}
\\
&=
\left[\mathrm{sc}^{i}(f)\right]_{i}
\tag{4}
\\
&=
\left[\mathrm{sc}^{i}(f)\right]_{\mathrm{nat}_{\mathsf{C}}(\mathrm{sc}^{i}(f))}
\tag{5}
\\
&=
\varepsilon^{[\mathsf{C}]}\left(\mathrm{sc}^{i}(f)\right).
\tag{6}
\end{align*}

The first equality unravels the definition of the mapping \(\varepsilon^{[\mathsf{C}]}\); the second equality follows from Proposition~\ref{PScTg}. Let us recall that from the explicit construction of the equivalence classes in \(\mathrm{Uf}^{(\omega)}(\mathrm{Gd}^{(\omega)}(\mathsf{C}))\), \([f]_{\mathrm{nat}_{\mathsf{C}}(f)}=[f]_{m}\). The third equality unravels the definition of the operation \(\mathrm{sc}'^{i}\) in \(\mathrm{Uf}^{(\omega)}(\mathrm{Gd}^{(\omega)}(\mathsf{C}))\); the fourth equality unravels the of the operation \(\mathrm{sc}^{(i,m)}\) in \(\mathrm{Gd}^{(\omega)}(\mathsf{C})\). Let us recall that \(\mathrm{sc}^{(i,m)}\) is just defined as the restriction \(\mathrm{sc}^{i}|_{C_{i}}^{C_{m}}\). The fifth equality follows from the explicit construction of the equivalence classes in \(\mathrm{Uf}^{(\omega)}(\mathrm{Gd}^{(\omega)}(\mathsf{C}))\); finally, the last equality recovers the definition of the mapping \(\varepsilon^{[\mathsf{C}]}\).

By a similar argument we also have that
\[
\varepsilon^{[\mathsf{C}]}\circ\mathrm{tg}'^{i}=\mathrm{tg}^{i}\circ\varepsilon^{[\mathsf{C}]}.
\]

\textsf{(Cii)}
Let \(f\) and \(g\) be elements of \(C\) with \(\mathrm{sc}^{i}(f)=\mathrm{tg}^{i}(g)\). Consider \(m=\max\{\mathrm{nat}_{\mathsf{C}(f)},\mathrm{nat}_{\mathsf{C}(g)},i\}+1\). The following chain of equalities holds
\begin{align*}
\varepsilon^{[\mathsf{C}]}\left(f\right)\#'^{i}\varepsilon^{[\mathsf{C}]}\left(g\right)
&=
\left[f\right]_{\mathrm{nat}_{\mathsf{C}}(f)}\#'^{i}\left[g\right]_{\mathrm{nat}_{\mathsf{C}}(g)}
\tag{1}
\\
&=
\left[f\right]_{m}\#'^{i}\left[g\right]_{m}
\tag{2}
\\
&=
\left[f\circ^{(i)}g\right]_{m}
\tag{3}
\\
&=
\left[f\#^{i}g\right]_{m}
\tag{4}
\\
&=
\left[f\#^{i}g\right]_{\mathrm{nat}_{\mathsf{C}}(f\#^{i}g)}
\tag{5}
\\
&=
\varepsilon^{[\mathsf{C}]}\left(f\#^{i}g\right).
\tag{6}
\end{align*}

The first equality unravels the definition of the mapping \(\varepsilon^{[\mathsf{C}]}\); the second equality follows from Proposition~\ref{PComp}. Let us recall that from the explicit construction of the equivalence classes in \(\mathrm{Uf}^{(\omega)}(\mathrm{Gd}^{(\omega)}(\mathsf{C}))\), \([f]_{\mathrm{nat}_{\mathsf{C}}(f)}=[f]_{m}\) and \([g]_{\mathrm{nat}_{\mathsf{C}}(g)}=[g]_{m}\). The third equality unravels the definition of the operation \(\#'^{i}\) in \(\mathrm{Uf}^{(\omega)}(\mathrm{Gd}^{(\omega)}(\mathsf{C}))\); the fourth equality unravels the definition of the operation \(\circ^{(i)}\) in \(\mathrm{Gd}^{(\omega)}(\mathsf{C})\). Let us recall that \(\circ^{(i)}\) is just defined as the restriction \(\#^{i}|_{C_{m}\times C_{m}}^{C_{m}}\). The fifth equality follows from the explicit construction of the equivalence classes in \(\mathrm{Uf}^{(\omega)}(\mathrm{Gd}^{(\omega)}(\mathsf{C}))\); finally, the last equality recovers the definition of the mapping \(\varepsilon^{[\mathsf{C}]}\).

All in all, \(\varepsilon^{[\mathsf{C}]}\) is a single-sorted \(\omega\)-functor from \(\mathsf{C}\) to \(\mathrm{Uf}^{(\omega)}(\mathrm{Gd}^{(\omega)}(\mathsf{C}))\),
\[
\varepsilon^{[\mathsf{C}]}:\mathsf{C}\mor\mathrm{Uf}^{(\omega)}(\mathrm{Gd}^{(\omega)}(\mathsf{C})).
\]

Claim~\ref{CUfwGdwB} follows.

Our next aim is to prove that, for every single-sorted \(\omega\)-category \(\mathsf{C}\), \(\eta^{[\mathsf{C}]}\) is a single-sorted \(\omega\)-functor.

\begin{claim}\label{CUfwGdwC}
\(\eta^{[\mathsf{C}]}\) is a single-sorted \(\omega\)-functor from \(\mathrm{Uf}^{(\omega)}(\mathrm{Gd}^{(\omega)}(\mathsf{C}))\) to \(\mathsf{C}\).
\end{claim}

It suffices to prove that \(\eta^{[\mathsf{C}]}\) satisfies the conditions of Definition~\ref{DnCatSS}.

Let \(i\in\omega\). We want to prove that \(\eta^{[\mathsf{C}]}\) is a single-sorted functor from \(([C_{\omega}], (\#'^{i},\mathrm{sc}'^{i},\mathrm{tg}'^{i}))\) to \((C, (\#^{i},\mathrm{sc}^{i},\mathrm{tg}^{i}))\). To this end, we consider the different items stated in Definition~\ref{DCat}.

\textsf{(Ci)}
Let \([f]_{j}\) be a \(\omega\)-cell class in \([C_{\omega}]\). Consider \(m=\max\{i,j\}+1\). The following chain of equalities holds
\begin{align*}
\eta^{[\mathsf{C}]}\left(\mathrm{sc}'^{i}\left(\left[f\right]_{j}\right)\right)
&=
\eta^{[\mathsf{C}]}\left(\mathrm{sc}'^{i}\left(\left[f\right]_{m}\right)\right)
\tag{1}
\\
&=
\eta^{[\mathsf{C}]}\left(\left[\mathrm{sc}^{(i,m)}(f)\right]_{i}\right)
\tag{2}
\\
&=
\mathrm{sc}^{(i,m)}(f)
\tag{3}
\\
&=
\mathrm{sc}^{i}(f)
\tag{4}
\\
&=
\mathrm{sc}^{i}\left(\eta^{[\mathsf{C}]}\left(\left[f\right]_{j}\right)\right).
\tag{5}
\end{align*}

The first equality follows from Proposition~\ref{PScTg}. Let us recall that from the explicit construction of the equivalence classes in \(\mathrm{Uf}^{(\omega)}(\mathrm{Gd}^{(\omega)}(\mathsf{C}))\), \([f]_{j}=[f]_{m}\). The second equality unravels the definition of \(\mathrm{sc}'^{i}\) in \(\mathrm{Uf}^{(\omega)}(\mathrm{Gd}^{(\omega)}(\mathsf{C}))\); the third equality unravels the definition of the mapping \(\eta^{[\mathsf{C}]}\); the fourth equality unravels the definition of \(\mathrm{sc}^{(i,m)}\) in \(\mathrm{Gd}^{(\omega)}(\mathsf{C})\). Let us recall that \(\mathrm{sc}^{(i,m)}\) is just defined as the restriction \(\mathrm{sc}^{i}|_{C_{m}}^{C_{i}}\). Finally, the last equality recovers the definition of the mapping \(\eta^{[\mathsf{C}]}\).

By a similar argument we also have that 
\[
\eta^{[\mathsf{C}]}\circ\mathrm{tg}'^{i}=\mathrm{tg}^{i}\circ\eta^{[\mathsf{C}]}.
\]

\textsf{(Cii)}
Let \([f]_{j}\) and \([g]_{k}\) be \(\omega\)-cell classes in \([C_{\omega}]\) with \(\mathrm{sc}'^{i}([f]_{j})=\mathrm{tg}'^{i}([g]_{k})\). Consider \(m=\max\{i,j,k\}+1\). The following chain of equalities holds
\allowdisplaybreaks
\begin{align*}
\eta^{[\mathsf{C}]}\left([f]_{j}\#'^{i}[g]_{k}\right)
&=
\eta^{[\mathsf{C}]}\left(\left[f\right]_{m}\#'^{i}\left[g\right]_{m}\right)
\tag{1}
\\
&=
\eta^{[\mathsf{C}]}\left(\left[f\circ^{(i)}g\right]_{m}\right)
\tag{2}
\\
&=
f\circ^{(i)}g
\tag{3}
\\
&=
f\#^{i}g
\tag{4}
\\
&=
\eta^{[\mathsf{C}]}\left([f]_{j}\right)\#^{i}\eta^{[\mathsf{C}]}\left([g]_{k}\right).
\tag{5}
\end{align*}

The first equality follows from Proposition~\ref{PComp}. Let us recall that from the explicit construction of the equivalence classes in \(\mathrm{Uf}^{(\omega)}(\mathrm{Gd}^{(\omega)}(\mathsf{C}))\), \([f]_{j}=[f]_{m}\) and \([g]_{k}=[g]_{m}\). The second equality unravels the definition of \(\#'^{i}\) in \(\mathrm{Uf}^{(\omega)}(\mathrm{Gd}^{(\omega)}(\mathsf{C}))\); the third equality unravels the definition of the mapping \(\eta^{[\mathsf{C}]}\); the fourth equality unravels the definition of \(\#^{(i)}\) in \(\mathrm{Gd}^{(\omega)}(\mathsf{C})\). Let us recall that \(\#^{(i)}\) is just defined as the restriction \(\#^{i}|_{C_{m}\times C_{m}}^{C_{i}}\). Finally, the last equality recovers the definition of the mapping \(\eta^{[\mathsf{C}]}\).

All in all, \(\eta^{[\mathsf{C}]}\) is a single-sorted \(\omega\)-functor from \(\mathrm{Uf}^{(\omega)}(\mathrm{Gd}^{(\omega)}(\mathsf{C}))\) to \(\mathsf{C}\),
\[
\eta^{[\mathsf{C}]}:\mathrm{Uf}^{(\omega)}(\mathrm{Gd}^{(\omega)}(\mathsf{C}))\mor\mathsf{C}.
\]

Claim~\ref{CUfwGdwC} follows.

We next prove that, for every single-sorted \(\omega\)-category \(\mathsf{C}\) in \(\w\mathsf{Cat}\), the single-sorted \(\omega\)-functor \(\varepsilon^{[\mathsf{C}]}\) and \(\eta^{[\mathsf{C}]}\) are mutually inverse morphisms in \(\w\mathsf{Cat}\).

\begin{claim}\label{CUfwGdwD}
For every single-sorted \(\omega\)-category \(\mathsf{C}\) in \(\w\mathsf{Cat}\),
\begin{align*}
\eta^{[\mathsf{C}]}\circ\varepsilon^{[\mathsf{C}]}&=\mathrm{Id}^{\mathsf{C}}
&\varepsilon^{[\mathsf{C}]}\circ\eta^{[\mathsf{C}]}&=\mathrm{Id}^{\mathrm{Uf}^{(\omega)}(\mathrm{Gd}^{(\omega)}(\mathsf{C}))}.
\end{align*}
\end{claim}

For every \(f\in C\), the following chain of equalities holds
\begin{align*}
\eta^{[\mathsf{C}]}\left(\varepsilon^{[\mathsf{C}]}\left(f\right)\right)
&=
\eta^{[\mathsf{C}]}\left(\left[f\right]_{\mathrm{nat}_{\mathsf{C}}(f)}\right)
\tag{1}
\\
&=
f.
\tag{2}
\end{align*}

The first equality unravels the definition of the single-sorted \(\omega\)-functor \(\varepsilon^{[\mathsf{C}]}\); the second equality unravels the definition of the single-sorted \(\omega\)-functor \(\eta^{[\mathsf{C}]}\).

Thus, \(\eta^{[\mathsf{C}]}\circ\varepsilon^{[\mathsf{C}]}=\mathrm{Id}^{\mathsf{C}}\).

For every \([f]_{j}\in[C_{\omega}]\), the following chain of equalities holds
\begin{align*}
\varepsilon^{[\mathsf{C}]}\left(\eta^{[\mathsf{C}]}\left(\left[f\right]_{j}\right)\right)
&=
\varepsilon^{[\mathsf{C}]}\left(f\right)
\tag{1}
\\
&=
[f]_{\mathrm{nat}_{\mathsf{C}}(f)}
\tag{2}
\\
&=
[f]_{j}.
\tag{3}
\end{align*}

The first equality unravels the definition of the single-sorted \(\omega\)-functor \(\eta^{[\mathsf{C}]}\); the second equality unravels the definition of the single-sorted \(\omega\)-functor \(\varepsilon^{[\mathsf{C}]}\); finally, the last equality follows from the explicit construction of the equivalence classes in \(\mathrm{Uf}^{(\omega)}(\mathrm{Gd}^{(\omega)}(\mathsf{C}))\).

Thus, \(\varepsilon^{[\mathsf{C}]}\circ\eta^{[\mathsf{C}]}=\mathrm{Id}^{\mathrm{Uf}^{(\omega)}(\mathrm{Gd}^{(\omega)}(\mathsf{C}))}\).

Claim~\ref{CUfwGdwD} follows.

We next prove that the morphisms \(\varepsilon^{[\bigcdot]}=(\varepsilon^{[\mathsf{C}]})_{\mathsf{C}\in\mathrm{Ob}(\w\mathsf{Cat})}\), when \(\mathsf{C}\) ranges over the objects of \(\w\mathsf{Cat}\), are the components of a natural transformation.

\begin{claim}\label{CUfwGdwE}
For every single-sorted \(\omega\)-functor \(F:\mathsf{C}\mor\mathsf{D}\), the square of single-sorted \(\omega\)-functors in Figure~\ref{FEtaNatTrans} commutes.
\begin{figure}[t]
\[
\xymatrix{
\mathsf{C}
  \ar[d]_{F}
  \ar[r]^-{\varepsilon^{[\mathsf{C}]}}
&
\mathrm{Uf}^{(\omega)}(\mathrm{Gd}^{(\omega)}(\mathsf{C}))
  \ar[d]^{\mathrm{Uf}^{(\omega)}(\mathrm{Gd}^{(\omega)}(F))}
\\
\mathsf{D}
  \ar[r]_-{\varepsilon^{[\mathsf{D}]}}
&
\mathrm{Uf}^{(\omega)}(\mathrm{Gd}^{(\omega)}(\mathsf{D}))
}
\]
\caption{The natural transformation \(\eta\) on single-sorted \(\omega\)-functors.}
\label{FEtaNatTrans}
\end{figure}
\end{claim}

Let \(\mathsf{C}\) and \(\mathsf{D}\) be the single-sorted \(\omega\)-categories under consideration, where
\begin{align*}
\mathsf{C}&=
\left(
C,
\left(
\#^{i\mathsf{C}},\mathrm{sc}^{i\mathsf{C}},\mathrm{tg}^{i\mathsf{C}}
\right)_{i\in\omega}
\right) \mbox{ and}
\\
\mathsf{D}&=
\left(
D,
\left(
\#^{i\mathsf{D}},\mathrm{sc}^{i\mathsf{D}},\mathrm{tg}^{i\mathsf{D}}
\right)_{i\in\omega}
\right),
\end{align*}
and let
\begin{align*}
\mathrm{Uf}^{(\omega)}\left(\mathrm{Gd}^{(\omega)}\left(\mathsf{C}\right)\right)&=
\left(
[C_{\omega}],
\left(
\#'^{i\mathsf{C}},\mathrm{sc}'^{i\mathsf{C}},\mathrm{tg}'^{i\mathsf{C}}
\right)_{i\in\omega}
\right) \mbox{ and}
\\
\mathrm{Uf}^{(\omega)}\left(\mathrm{Gd}^{(\omega)}\left(\mathsf{D}\right)\right)&=
\left(
[D_{\omega}],
\left(
\#'^{i\mathsf{D}},\mathrm{sc}'^{i\mathsf{D}},\mathrm{tg}'^{i\mathsf{D}}
\right)_{i\in\omega}
\right)
\end{align*}
be the corresponding unifications of the graduations of \(\mathsf{C}\) and \(\mathsf{D}\), respectively.

Let \(F\) be a single-sorted \(\omega\)-functor from \(\mathsf{C}\) to \(\mathsf{D}\). Let us recall that, according to Remark~\ref{RGdw}, its graduation is given by \(\mathrm{Gd}^{(\omega)}(F)=(F|^{D_{k}}_{C_{k}})_{k\in\omega}\). Moreover, according to Definition~\ref{DUfwFun}, the unification of the graduation of \(F\) is the mapping that assigns to every \(\omega\)-cell class \([f]_{j}\) in \([C_{\omega}]\) the \(\omega\)-cell class \(\left[F|^{D_{j}}_{C_{j}}(f)\right]_{j}\) in \([D_{\omega}]\).

Let \(\varepsilon^{[\mathsf{C}]}\) be the single-sorted \(\omega\)-functor from \(\mathsf{C}\) to \(\mathrm{Uf}^{(\omega)}(\mathrm{Gd}^{(\omega)}(\mathsf{C}))\) and \(\varepsilon^{[\mathsf{D}]}\) the single-sorted \(\omega\)-functor from \(\mathsf{D}\) to \(\mathrm{Uf}^{(\omega)}(\mathrm{Gd}^{(\omega)}(\mathsf{D}))\) defined before. Let us recall that, for every \(f\in C\) and every \(g\in D\), \(\varepsilon^{[\mathsf{C}]}\) and \(\varepsilon^{[\mathsf{D}]}\) are defined as
\begin{align*}
\varepsilon^{[\mathsf{C}]}(f)&=[f]_{\mathrm{nat}_{\mathsf{C}}(f)};
&\varepsilon^{[\mathsf{D}]}(g)&=[g]_{\mathrm{nat}_{\mathsf{D}}(g)}.
\end{align*}

Let \(f\in C\), then the following chain of equalities holds
\begin{align*}
\mathrm{Uf}^{(\omega)}(\mathrm{Gd}^{(\omega)}(F))\left(\varepsilon^{[\mathsf{C}]}(f)\right)
&=
\mathrm{Uf}^{(\omega)}(\mathrm{Gd}^{(\omega)}(F))\left(\left[f\right]_{\mathrm{nat}_{\mathsf{C}}(f)}\right)
\tag{1}
\\
&=
\left[F\big|^{D_{\mathrm{nat}_{\mathsf{C}}(f)}}_{C_{\mathrm{nat}_{\mathsf{C}}(f)}}(f)\right]_{\mathrm{nat}_{\mathsf{C}}(f)}
\tag{2}
\\
&=
\left[F(f)\right]_{\mathrm{nat}_{\mathsf{C}}(f)}
\tag{3}
\\
&=
\left[F(f)\right]_{\mathrm{nat}_{\mathsf{D}}(F(f))}
\tag{4}
\\
&=
\varepsilon^{[\mathsf{D}]}\left(F(f)\right).
\tag{5}
\end{align*}

The first equality unravels the definition of the single-sorted \(\omega\)-functor \(\varepsilon^{[\mathsf{C}]}\); the second equality unravels the definition of the single-sorted \(\omega\)-functor \(\mathrm{Uf}^{(\omega)}\circ\mathrm{Gd}^{(\omega)}(F)\); the third equality applies the definition of a restriction of a mapping; the fourth equality follows from the fact that, since \(F\) is a single-sorted \(\omega\)-functor, then \(\mathrm{nat}_{\mathsf{C}}(f)=\mathrm{nat}_{\mathsf{D}}(F(f))\) by item~(MSii) in Definition~\ref{DnCatMS}; finally, the last equality recovers the definition of the single-sorted \(\omega\)-functor \(\varepsilon^{[\mathsf{D}]}\).

Claim~\ref{CUfwGdwE} follows.

So, considering the foregoing, we can affirm that \(\varepsilon\) is a natural transformation from \(\mathrm{Id}^{\w\mathsf{Cat}}\) to \(\mathrm{Gd}^{(\omega)}\circ\mathrm{Uf}^{(\omega)}\).

Let us note that each single-sorted \(\omega\)-category \(\mathsf{C}\) in \(\w\mathsf{Cat}\) determines, by Claim~\ref{CUfwGdwB}, the single-sorted \(\omega\)-functor \(\varepsilon^{[\mathsf{C}]}\) which, by Claim~\ref{CUfwGdwD}, is an isomorphism.

Thus, \(\varepsilon^{[\bigcdot]}:\mathrm{Id}^{\w\mathsf{Cat}}\cel\mathrm{Uf}^{(\omega)}\circ\mathrm{Gd}^{(\omega)}\) is a natural isomorphism.
\end{proof}

Moreover, similar to the previous result, the process of \(\omega\)-unification followed by that of \(\omega\)-graduation is naturally isomorphic to the identity functor on the category of many-sorted \(\omega\)-categories.

\begin{proposition}\label{PGdwUfw} 
$\mathrm{Gd}^{(\omega)}\circ\mathrm{Uf}^{(\omega)}\cong\mathrm{Id}^{\w\mathsf{Cat}^{\mathrm{MS}}}$.
\end{proposition}

\begin{proof}
Let \(\mathsf{C}=((C_{k})_{k\in\omega},\zeta)\) be a many-sorted \(\omega\)-category with \(\zeta=(\zeta_{j,k})_{(j,k)\in\coprod_{k\in\omega}k}\) and, for every \((j,k)\in\coprod_{k\in\omega}k\),
\[
\zeta_{j,k}=
\left(
\circ^{(j)},\mathrm{sc}^{(j,k)},\mathrm{tg}^{(j,k)},\mathrm{i}^{(k,j)}
\right),
\]
and \(\mathrm{Uf}^{(\omega)}(\mathsf{C})=([C_{\omega}],(\#^{i},\mathrm{sc}^{i},\mathrm{tg}^{i})_{i\in\omega})\) its unification where \([C_{\omega}]\) is the quotient set \(C_{\omega}/\Theta\) presented in Definition~\ref{D[Cw]} and its structural operations are defined by cases in Definition~\ref{DScTg} and \ref{DComp}.

Now, let us consider the many-sorted \(\omega\)-category given by the graduation of the unification of \(\mathsf{C}\), i.e.,
\[
\mathrm{Gd}^{(\omega)}\left(\mathrm{Uf}^{(\omega)}\left(\mathsf{C}\right)\right)=
\left(
\left(\left[C_{\omega}\right]_{k}\right)_{k\in\omega},
\left(
\circ'^{(j)},\mathrm{sc}'^{(j,k)},\mathrm{tg}'^{(j,k)},\mathrm{i}'^{(k,j)}
\right)_{(j,k)\in\coprod_{k\in\omega}k}
\right)
\]
where according to Remark~\ref{RGdw}, for every \(k\in\omega\),
\[
[C_{\omega}]_{k}=\left\{[f]_{j}\in [C_{\omega}]\mid [f]_{j}=\mathrm{sc}^{k}([f]_{j})\right\}.
\]

Moreover, for every \(0\leq j<k\), the structural operations of the graduation are given by
\begin{align*}
\circ'^{(j)}&=\#^{j}\big|^{[C_{\omega}]_{k}}_{[C_{\omega}]_{k}\times [C_{\omega}]_{k}};
&\mathrm{i}'^{(k,j)}&=\mathrm{in}^{[C_{\omega}]_{j},[C_{\omega}]_{k}};
\\
\mathrm{sc}'^{(j,k)}&=\mathrm{sc}^{j}\big|^{[C_{\omega}]_{j}}_{[C_{\omega}]_{k}};
&\mathrm{tg}'^{(j,k)}&=\mathrm{tg}^{j}\big|^{[C_{\omega}]_{j}}_{[C_{\omega}]_{k}}.
\end{align*}

We let \(\alpha^{[\mathsf{C}]}=\left(\alpha_{i}^{[\mathsf{C}]}\right)_{i\in\omega}\) stand for the \(\omega\)-sorted mapping where, for every \(i\in\omega\), \(\alpha^{[\mathsf{C}]}_{i}\) is the mapping \(\mathrm{in}^{[C_{i}]}=\mathrm{pr}^{\Theta}\circ\mathrm{in}^{C_{i}}\) from \(C_{i}\) to \([C_{\omega}]\) presented in Definition~\ref{D[Cw]}. Let us recall that \(\mathrm{in}^{[C_{i}]}\) sends an element \(f\in C_{i}\) to the \(\omega\)-cell class \([f]_{i}\).
\[
\alpha^{[\mathsf{C}]}_{i}:C_{i}\mor[C_{\omega}]
\]

Moreover, we let \(\beta^{\mathsf{C}}=(\beta^{\mathsf{C}}_{i})_{i\in\omega}\) stand for the \(\omega\)-sorted mapping where, for every \(i\in\omega\), \(\beta^{\mathsf{C}}_{i}\) is the assignment from \(C_{\omega}\) to \(C_{i}\) defined, for every \(\omega\)-cell \((f,j)\) in \(C_{\omega}\), as 
\[
\beta^{\mathsf{C}}_{i}(f,j)
=
\begin{cases}
\mathrm{i}^{(i,j)}(f) &\mbox{if \(j<i\);}
\\
f &\mbox{if \(j=i\);}
\\
\mathrm{sc}^{(i,j)}(f) &\mbox{if \(i<j\).}
\end{cases}
\]

By symmetry in the definition of many-sorted \(\omega\)-category, we could also define the mapping \(\beta^{\mathsf{C}}_{i}\) using the \((i,j)\)-target instead of the \((i,j)\)-source.

\begin{claim}\label{CGdwUfwA}
For every \(i\in\omega\) and every \(\omega\)-cells \((f,j)\) and \((g,k)\) in \(C_{\omega}\), if \(((f,j),(g,k))\in\Theta\), then \(\beta^{\mathsf{C}}_{i}(f,j)=\beta^{\mathsf{C}}_{i}(g,k)\).
\end{claim}

Let us first note that \(\mathrm{sc}^{i}=\mathrm{in}^{[C_{i}]}\circ\beta^{\mathsf{C}}_{i}\circ\mathrm{pr}^{\Theta}\). That is, for every \(\omega\)-cell \((f,j)\) in \(C_{\omega}\),
\[
\mathrm{sc}^{i}([f]_{j})=\left[\beta^{\mathsf{C}}_{i}(f,j)\right]_{i}
\]

If \(((f,j),(g,k))\in\Theta\), then, applying Proposition~\ref{PScTg},
\[
\left[\beta^{\mathsf{C}}_{i}(f,j)\right]_{i}
=
\mathrm{sc}^{i}([f]_{j})
=
\mathrm{sc}^{i}([g]_{k})
=
\left[\beta^{\mathsf{C}}_{i}(g,k)\right]_{i}.
\]
Therefore, according to the definition of \(\Theta\) presented in Definition~\ref{D[Cw]}, \(\beta^{\mathsf{C}}_{i}(f,j)=\beta^{\mathsf{C}}_{i}(g,k)\).

Claim~\ref{CGdwUfwA} follows.

We let \(\beta^{[\mathsf{C}]}_{i}=\left(\beta^{\mathsf{C}}_{i}\right)^{\natural}\) stand for the mapping obtained by the universal property of the quotient as the unique mapping from \([C_{\omega}]\) to \(C_{i}\) such that \(\left(\beta^{\mathsf{C}}_{i}\right)^{\natural}\circ\mathrm{pr}^{\Theta}=\beta^{\mathsf{C}}_{i}\). Thus, \(\beta^{[\mathsf{C}]}_{i}\) is the mapping from \([C_{\omega}]\) to \(C_{i}\) that assigns to a \(\omega\)-cell class \([f]_{j}\) in \([C_{\omega}]\) the element \(\eta^{\mathsf{C}}_{i}(f,j)\), i.e., the value of the mapping \(\eta^{\mathsf{C}}_{i}\) at any class representative. That \(\eta^{[\mathsf{C}]}_{i}\) is well-defined follows from Claim~\ref{CGdwUfwA}. 

\begin{claim}\label{CGdwUfwB}
For every \(i\in\omega\), \(\beta^{[\mathsf{C}]}_{i}\left[\left[C_{\omega}\right]_{i}\right]=C_{i}\).
\end{claim}

That \(\beta^{[\mathsf{C}]}_{i}\left[\left[C_{\omega}\right]_{i}\right]\subseteq C_{i}\) follows by definition of \(\beta^{\mathsf{C}}_{i}\). 

Conversely, let \(f\in C_{i}\). We want to show that \(f\in \beta^{[C]}_{i}\left[\left[C_{\omega}\right]_{i}\right]\).

Consider the \(\omega\)-cell class \([f]_{i}\). By definition of \(\mathrm{sc}^{i}\) presented in Definition~\ref{DScTg}, The equality \([f]_{i}=\mathrm{sc}^{i}\left([f]_{i}\right)\) follows. Therefore, \([f]_{i}\in[C_{\omega}]_{i}\). Moreover, unravelling the corresponding definitions
\[
\beta^{[\mathsf{C}]}_{i}\left([f]_{i}\right)
=
\beta^{\mathsf{C}}_{i}\left(f,i\right)
=
f.
\]

Claim~\ref{CGdwUfwB} follows.

Finally, we let \(\beta^{[\mathsf{C}]}\) stand for the \(\omega\)-sorted mapping \(\left(\beta^{[\mathsf{C}]}_{i}\right)_{i\in\omega}\).

Our next aim is to prove that, for every many-sorted \(\omega\)-category \(\mathsf{C}\), \(\alpha^{[\mathsf{C}]}\) is a many-sorted \(\omega\)-functor.

\begin{claim}\label{CGdwUfwC}
\(\alpha^{[\mathsf{C}]}\) is a many-sorted \(\omega\)-functor from \(\mathsf{C}\) to \(\mathrm{Gd}^{(\omega)}(\mathrm{Uf}^{(\omega)}(\mathsf{C}))\).
\end{claim}

It suffices to prove that \(\varepsilon^{\mathsf{C}}\) satisfies the conditions of Definition~\ref{DnCatMS}.

\textsf{(MSi)}
For every \(0\leq j<k\), if \(f\in C_{k}\), the following chain of equalities holds
\allowdisplaybreaks
\begin{align*}
\mathrm{sc}'^{(j,k)}\left(\alpha_{k}^{[\mathsf{C}]}(f)\right)
&=
\mathrm{sc}'^{(j,k)}\left(\left[f\right]_{k}\right)
\tag{1}
\\
&=
\mathrm{sc}^{j}\left(\left[f\right]_{k}\right)
\tag{2}
\\
&=
\left[\mathrm{sc}^{(j,k)}(f)\right]_{j}
\tag{3}
\\
&=
\alpha^{[\mathsf{C}]}_{j}\left(\mathrm{sc}^{(j,k)}(f)\right).
\tag{4}
\end{align*}

The first equality unravels the \(k\)-th component of the \(\omega\)-sorted mapping \(\alpha^{[\mathsf{C}]}\); the second equality unravels the definition of the structural operation \(\mathrm{sc}'^{(j,k)}\) in \(\mathrm{Gd}^{(\omega)}(\mathrm{Uf}^{(\omega)}(\mathsf{C}))\). Let us recall that \(\mathrm{sc}'^{(j,k)}\) is just defined as the birestriction \(\mathrm{sc}^{j}|^{[C_{\omega}]_{j}}_{[C_{\omega}]_{k}}\). The third equality unravels the definition of the structural operation \(\mathrm{sc}^{j}\) in \(\mathrm{Uf}^{(\omega)}(\mathsf{C})\) presented in Definition~\ref{DScTg}; finally, the last equality recovers the definition of the \(j\)-th component of the \(\omega\)-sorted mapping \(\alpha^{[\mathsf{C}]}\).

By a similar argument we also have that 
\[
\mathrm{tg}'^{(j,k)}\circ\alpha^{[\mathsf{C}]}_{k}
=
\alpha^{[\mathsf{C}]}_{j}\circ\mathrm{tg}^{(j,k)}.
\]

\textsf{(MSii)}
For every \(0\leq j<k\), if \(f\in C_{j}\), the following chain of equalities holds
\allowdisplaybreaks
\begin{align*}
\mathrm{i}'^{(k,j)}\left(\alpha^{[\mathsf{C}]}_{j}(f)\right)
&=
\mathrm{i}'^{(k,j)}\left(\left[f\right]_{j}\right)
\tag{1}
\\
&=
[f]_{j}
\tag{2}
\\
&=
\left[\mathrm{i}^{(k,l)}(f)\right]_{k}
\tag{3}
\\
&=
\alpha^{[\mathsf{C}]}_{k}\left(\mathrm{i}^{(k,j)}\left(f\right)\right).
\tag{4}
\end{align*}

The first equality unravels the definition of the \(j\)-th component of the \(\omega\)-sorted mapping \(\alpha^{[\mathsf{C}]}\); the second equality unravels the definition of the structural operation \(\mathrm{i}'^{(k,j)}\) in \(\mathrm{Gd}^{(\omega)}(\mathrm{Uf}^{(\omega)}(\mathsf{C}))\). Let us recall that \(\mathrm{i}'^{(k,j)}\) is just defined as the inclusion \(\mathrm{in}^{[C_{\omega}]_{j},[C_{\omega}]_{k}}\). The third equality follows from the definition of the equivalence relation \(\Theta\) presented in Definition~\ref{D[Cw]}; finally, the last equality recovers the definition of the \(k\)-th component of the \(\omega\)-sorted mapping \(\alpha^{[\mathsf{C}]}\).

\textsf{(MSiii)}
For every \(0\leq j<k\), if \(f,g\in C_{k}\) with \(\mathrm{sc}^{(j,k)}(f)=\mathrm{tg}^{(j,k)}(g)\), the following chain of equalities holds
\begin{align*}
\alpha^{[\mathsf{C}]}_{k}(f)\circ'^{(j)}\alpha^{[\mathsf{C}]}_{k}(g)
&=
[f]_{k}\circ'^{(j)}[g]_{k}
\tag{1}
\\
&=
[f]_{k}\#^{i}[g]_{k}
\tag{2}
\\
&=
[f\circ^{(j)}g]_{k}
\tag{3}
\\
&=
\alpha^{[\mathsf{C}]}_{k}(f\circ^{(j)}g).
\tag{4}
\end{align*}

The first equality unravels the definition of the \(k\)-th component of the \(\omega\)-sorted mapping \(\alpha^{[\mathsf{C}]}\); the second equality unravels the definition of the structural operation \(\circ'^{(j)}\) in \(\mathrm{Gd}^{(\omega)}(\mathrm{Uf}^{(\omega)}(\mathsf{C}))\). Let us recall that \(\circ'^{(j)}\) is just defined as the birestriction \(\#^{j}|_{[C_{\omega}]_{k}\times[C_{\omega}]_{k}}^{[C_{\omega}]_{k}}\). The third equality unravels the definition of the structural operation \(\#^{i}\) in \(\mathrm{Uf}^{(\omega)}(\mathsf{C})\) presented in Definition~\ref{DComp}; finally, the last equality recovers the definition of the \(k\)-th component of the \(\omega\)-sorted mapping \(\alpha^{[\mathsf{C}]}\).

All in all, \(\alpha^{[\mathsf{C}]}\) is a many-sorted \(\omega\)-functor from \(\mathsf{C}\) to \(\mathrm{Gd}^{(\omega)}(\mathrm{Uf}^{(\omega)}(\mathsf{C}))\),
\[
\alpha^{[\mathsf{C}]}:\mathsf{C}\mor\mathrm{Gd}^{(\omega)}(\mathrm{Uf}^{(\omega)}(\mathsf{C})).
\]

Claim~\ref{CGdwUfwC} follows.

Before proceeding any further, let us point out the following fact. Let \(k\in\omega\) and \([f]_{j}\in[C_{\omega}]\) in its natural representation, in particular \(j=\mathrm{nat}_{\mathsf{C}}([f]_{j})\) is the nature of \([f]_{j}\) in \(\mathsf{C}\). If \([f]_{j}\in[C_{\omega}]_{k}\), then \(j\leq k\). If \(k<j\), then the following chain of equalities holds
\[
[f]_{j}
=
\mathrm{sc}^{k}\left([f]_{j}\right)
=
\left[\mathrm{sc}^{(k,j)}(f)\right]_{k},
\]
in contradiction with \(j\) being the nature of \([f]_{j}\).

Therefore, for every \(\omega\)-cell class \([f]_{j}\) in \([C_{\omega}]\), if \([f]_{j}\in[C_{\omega}]_{k}\), i.e., if \([f]_{j}\) is a \(k\)-cell, that is, \([f]_{j}=\mathrm{sc}^{k}\left([f]_{j}\right)\), then \(\mathrm{nat}_{\mathsf{C}}([f]_{j})<k\). Thus, we can affirm that
\[
[C_{\omega}]_{k}=\{[f]_{j}\in [C_{\omega}]\mid \mathrm{nat}_{\mathsf{C}}([f]_{j})<k\}.
\]

Our next aim is to prove that, for every many-sorted \(\omega\)-category \(\mathsf{C}\), \(\beta^{[\mathsf{C}]}\) is a many-sorted \(\omega\)-functor.

\begin{claim}\label{CGdwUfwD}
\(\beta^{[\mathsf{C}]}\) is a many-sorted \(\omega\)-functor from \(\mathrm{Gd}^{(\omega)}(\mathrm{Uf}^{(\omega)}(\mathsf{C}))\) to \(\mathsf{C}\).
\end{claim}

It suffices to prove that \(\beta^{[\mathsf{C}]}\) satisfies the conditions of Definition~\ref{DnCatMS}.

\textsf{(MSi)}
Let \(0\leq j<k\) and \([f]_{i}\in[C_{\omega}]_{k}\) in its natural representation. Therefore, \(i\leq k\). If \(i<k\), then the following chain of equalities holds
\begin{align*}
\beta^{[\mathsf{C}]}_{j}\left(\mathrm{sc}'^{(j,k)}\left(\left[f\right]_{i}\right)\right)
&=
\beta^{[\mathsf{C}]}_{j}\left(\mathrm{sc}^{j}\left(\left[f\right]_{i}\right)\right)
\tag{1}
\\
&=
\beta^{[\mathsf{C}]}_{j}\left(\mathrm{sc}^{j}\left(\left[\mathrm{i}^{(k,i)}(f)\right]_{k}\right)\right)
\tag{2}
\\
&=
\beta^{[\mathsf{C}]}_{j}\left(\left[\mathrm{sc}^{(j,k)}\left(\mathrm{i}^{(k,i)}(f)\right)\right]_{j}\right)
\tag{3}
\\
&=
\mathrm{sc}^{(j,k)}\left(\mathrm{i}^{(k,i)}(f)\right)
\tag{4}
\\
&=
\mathrm{sc}^{(j,k)}\left(\beta^{[\mathsf{C}]}_{k}\left(\left[f\right]_{i}\right)\right).
\tag{5}
\end{align*}

The first equality unravels the definition of the structural operation \(\mathrm{sc}'^{(j,k)}\) in \(\mathrm{Gd}^{(\omega)}(\mathrm{Uf}^{(\omega)}(\mathsf{C}))\). Let us recall that \(\mathrm{sc}'^{(j,k)}\) is just defined as the birestriction \(\mathrm{sc}^{j}|^{[C_{\omega}]_{j}}_{[C_{\omega}]_{k}}\). The second equality follows from Proposition~\ref{PScTg}; the third equality unravels the definition of the structural operation \(\mathrm{sc}	^{j}\) in \(\mathrm{Uf}^{(\omega)}(\mathsf{C})\); The fourth equality unravels the definition of the \(j\)-th component of the \(\omega\)-sorted mapping \(\beta^{[\mathsf{C}]}\); finally, the last equality recovers the definition of the \(k\)-th component of the \(\omega\)-sorted mapping \(\beta^{[\mathsf{C}]}\).

If \(i=k\), the following chain of equalities holds
\begin{align*}
\beta^{[\mathsf{C}]}_{j}\left(\mathrm{sc}'^{(j,i)}\left(\left[f\right]_{i}\right)\right)
&=
\beta^{[\mathsf{C}]}_{j}\left(\mathrm{sc}^{j}\left(\left[f\right]_{i}\right)\right)
\tag{1}
\\
&=
\beta^{[\mathsf{C}]}_{j}\left(\left[\mathrm{sc}^{(j,i)}(f)\right]_{j}\right)
\tag{2}
\\
&=
\mathrm{sc}^{(j,i)}(f)
\tag{3}
\\
&=
\mathrm{sc}^{(j,k)}\left(\beta^{[\mathsf{C}]}_{i}\left(\left[f\right]_{i}\right)\right).
\tag{4}
\end{align*}

The first equality unravels the definition of the structural operation \(\mathrm{sc}'^{(j,k)}\) in \(\mathrm{Gd}^{(\omega)}(\mathrm{Uf}^{(\omega)}(\mathsf{C}))\). Let us recall that \(\mathrm{sc}'^{(j,i)}\) is just defined as the birestriction \(\mathrm{sc}^{j}|^{[C_{\omega}]_{j}}_{[C_{\omega}]_{i}}\). The second equality unravels the definition of the structural operation \(\mathrm{sc}	^{j}\) in \(\mathrm{Uf}^{(\omega)}(\mathsf{C})\); the fourth equality unravels the definition of the \(j\)-th component of the \(\omega\)-sorted mapping \(\beta^{[\mathsf{C}]}\); finally, the last equality recovers the definition of the \(i\)-th component of the \(\omega\)-sorted mapping \(\beta^{[\mathsf{C}]}\).

By a similar argument we also have that 
\[
\mathrm{tg}^{(j,k)}\circ\beta^{[\mathsf{C}]}_{k}
=
\beta^{[\mathsf{C}]}_{j}\circ\mathrm{tg}'^{(j,k)}.
\]

\textsf{(MSii)}
Let \(0\leq j<k\) and \([f]_{i}\in[C_{\omega}]_{j}\) in its natural representation. Therefore, \(i\leq j<k\). If \(i<j\), the following chain of equalities holds
\begin{align*}
\mathrm{i}^{(k,j)}\left(\beta^{[\mathsf{C}]}_{j}\left([f]_{i}\right)\right)
&=
\mathrm{i}^{(k,j)}\left(\mathrm{i}^{(j,i)}(f)\right)
\tag{1}
\\
&=
\mathrm{i}^{(k,i)}(f)
\tag{2}
\\
&=
\beta^{[\mathsf{C}]}_{k}\left([f]_{i}\right)
\tag{3}
\\
&=
\beta^{[\mathsf{C}]}_{k}\left(\mathrm{i}'^{(k,j)}([f]_{i})\right).
\tag{4}
\end{align*}

The first equality unravels the definition of the \(j\)-th component of the \(\omega\)-sorted mapping \(\beta^{[\mathsf{C}]}\) since \(i<j\); the second equality follows from item~(MS4) in Definition~\ref{DnCatMS}; the third equality recovers the definition of the \(k\)-th component of the \(\omega\)-sorted mapping \(\beta^{[\mathsf{C}]}\) since \(i<k\); finally, the last equality recovers the definition of the structural operation \(\mathrm{i}'^{(k,j)}\) in \(\mathrm{Gd}^{(\omega)}(\mathrm{Uf}^{(\omega)}(\mathsf{C}))\). Let us recall that \(\mathrm{i}'^{(k,j)}\) is just defined as the inclusion \(\mathrm{in}^{[C_{\omega}]_{j},[C_{\omega}]_{k}}\).

Moreover, if \(j=i\), the following chain of equalities holds
\begin{align*}
\mathrm{i}^{(k,i)}\left(\beta^{[\mathsf{C}]}_{i}\left([f]_{i}\right)\right)
&=
\mathrm{i}^{(k,i)}\left(f\right)
\tag{1}
\\
&=
\beta^{[\mathsf{C}]}_{k}\left([f]_{i}\right)
\tag{3}
\\
&=
\beta^{[\mathsf{C}]}_{k}\left(\mathrm{i}'^{(k,i)}\left([f]_{i}\right)\right).
\tag{4}
\end{align*}

The first equality unravels the definition of the \(j\)-th component of the \(\omega\)-sorted mapping \(\beta^{[\mathsf{C}]}\), since \(i=j\); the second equality recovers the definition of the \(k\)-th component of the \(\omega\)-sorted mapping \(\beta^{[\mathsf{C}]}\), since \(i<k\); finally, the last equality recovers the definition of the structural operation \(\mathrm{i}'^{(k,j)}\) in \(\mathrm{Gd}^{(\omega)}(\mathrm{Uf}^{(\omega)}(\mathsf{C}))\). Let us recall that \(\mathrm{i}'^{(k,j)}\) is just defined as the inclusion \(\mathrm{in}^{[C_{\omega}]_{j},[C_{\omega}]_{k}}\).

\textsf{(MSiii)}
Let \(0\leq j<m\) and \([f]_{l},[g]_{m}\in [C_{\omega}]_{k}\) in their natural representation with \(\mathrm{sc}'^{(j,k)}([f]_{l})=\mathrm{tg}'^{(j,k)}([g]_{m})\). Therefore, \(l\leq k\) and \(m\leq k\).

If \(l<k\) and \(m<k\), the following chain of equalities holds
\begin{align*}
\beta^{[\mathsf{C}]}_{k}\left([f]_{l}\circ'^{(j)}[g]_{m}\right)
&=
\beta^{[\mathsf{C}]}_{k}\left([f]_{l}\#^{j}[g]_{m}\right)
\tag{1}
\\
&=
\beta^{[\mathsf{C}]}_{k}\left([\mathrm{i}^{(k,l)}(f)]_{k}\#^{j}[\mathrm{i}^{(k,m)}(g)]_{k}\right)
\tag{2}
\\
&=
\beta^{[\mathsf{C}]}_{k}\left([\mathrm{i}^{(k,l)}(f)\circ^{(j)}\mathrm{i}^{(k,m)}(g)]_{k}\right)
\tag{3}
\\
&=
\mathrm{i}^{(k,l)}(f)\circ^{(j)}\mathrm{i}^{(k,m)}(g)
\tag{4}
\\
&=
\beta^{[\mathsf{C}]}\left([f]_{l}\right)\circ^{(j)}\beta^{[\mathsf{C}]}\left([g]_{m}\right).
\tag{5}
\end{align*}

The first equality unravels the definition of the structural operation \(\circ'^{(j)}\) in \(\mathrm{Gd}^{(\omega)}(\mathrm{Uf}^{(\omega)}(\mathsf{C}))\). Let us recall that \(\circ'^{(j)}\) is just defined as the birestriction \(\#^{j}|^{[C_{\omega}]_{k}}_{[C_{\omega}]_{k}\times[C_{\omega}]_{k}}\). The second equality follows from Proposition~\ref{PComp}; the third equality unravels the definition of the structural operation \(\#^{j}\) in \(\mathrm{Uf}^{(\omega)}(\mathsf{C})\) presented in Definition~\ref{DComp}; the fourth equality unravels the definition of \(k\)-th component of the \(\omega\)-sorted mapping \(\beta^{[\mathsf{C}]}\); finally, the last equality recovers the definition of the \(k\)-th component of the \(\omega\)-sorted mapping \(\beta^{[\mathsf{C}]}\).

The rest of the cases follow in a similar way.

All in all, \(\beta^{[\mathsf{C}]}\) is a many-sorted \(\omega\)-functor from \(\mathrm{Gd}^{(\omega)}(\mathrm{Uf}^{(\omega)}(\mathsf{C}))\) to \(\mathsf{C}\),
\[
\beta^{[\mathsf{C}]}:\mathrm{Gd}^{(\omega)}(\mathrm{Uf}^{(\omega)}(\mathsf{C}))\mor\mathsf{C}.
\]

Claim~\ref{CGdwUfwD} follows.

We next prove that, for every many-sorted \(\omega\)-category \(\mathsf{C}\) in \(\w\mathsf{Cat}^{\mathrm{MS}}\), the many-sorted \(\omega\)-functors \(\alpha^{[\mathsf{C}]}\) and \(\beta^{[\mathsf{C}]}\) are mutually inverse morphisms in \(\w\mathsf{Cat}^{\mathrm{MS}}\).

\begin{claim}\label{CGdwUfwE}
For every many-sorted \(\omega\)-category \(\mathsf{C}\) in \(\w\mathsf{Cat}^{\mathrm{MS}}\),
\begin{align*}
\beta^{[\mathsf{C}]}\circ\alpha^{[\mathsf{C}]}&=\mathrm{Id}^{\mathsf{C}}
&\alpha^{[\mathsf{C}]}\circ\beta^{[\mathsf{C}]}&=\mathrm{Id}^{\mathrm{Gd}^{(\omega)}(\mathrm{Uf}^{(\omega)}(\mathsf{C}))}.
\end{align*}
\end{claim}

For every \(i\in\omega\), if \(f\in C_{i}\), the following chain of equalities holds
\begin{align*}
\beta^{[\mathsf{C}]}_{i}\left(\alpha^{[\mathsf{C}]}_{i}\left(f\right)\right)
&=
\beta^{[\mathsf{C}]}_{i}\left(\left[f\right]_{i}\right)
\tag{1}
\\
&=
f.
\tag{2}
\end{align*}

The first equality unravels the definition of the \(i\)-th component of \(\alpha^{[\mathsf{C}]}\); the second equality unravels the definition of the \(i\)-th component of \(\beta^{[\mathsf{C}]}\).

Thus, \(\beta^{[\mathsf{C}]}\circ\alpha^{[\mathsf{C}]}=\mathrm{Id}^{\mathsf{C}}\).

Let \(i\in\omega\) and \([f]_{j}\in[C_{\omega}]_{i}\) in its natural representation. Therefore, \(j\leq i\). 

If \(j<i\), the following chain of equalities holds
\begin{align*}
\alpha^{[\mathsf{C}]}_{i}\left(\beta^{[\mathsf{C}]}_{i}\left([f]_{j}\right)\right)
&=
\alpha^{[\mathsf{C}]}_{i}\left(\mathrm{i}^{(i,j)}(f)\right)
\tag{1}
\\
&=
\left[\mathrm{i}^{(i,j)}(f)\right]_{i}
\tag{2}
\\
&=
[f]_{j}.
\tag{3}
\end{align*}

The first equality unravels the \(i\)-th component of the many-sorted \(\omega\)-functor \(\beta^{[\mathsf{C}]}\); the second equality unravels the \(i\)-th component of the many-sorted \(\omega\)-functor \(\alpha^{[\mathsf{C}]}\); finally, the last equality follows from the definition of the equivalence relation \(\Theta\) presented	in Definition~\ref{D[Cw]}.

If \(j=i\), then the following chain of equalities holds
\begin{align*}
\alpha^{[\mathsf{C}]}_{i}\left(\beta^{[\mathsf{C}]}_{i}\left([f]_{i}\right)\right)
&=
\alpha^{[\mathsf{C}]}_{i}\left(f\right)
\tag{1}
\\
&=
[f]_{i}
\end{align*}

The first equality unravels the \(i\)-th component of the many-sorted \(\omega\)-functor \(\beta^{[\mathsf{C}]}\); the second equality unravels the \(i\)-th component of the many-sorted \(\omega\)-functor \(\alpha^{[\mathsf{C}]}\).

Thus, \(\alpha^{[\mathsf{C}]}\circ\beta^{[\mathsf{C}]}=\mathrm{Id}^{\mathrm{Gd}^{(\omega)}(\mathrm{Uf}^{(\omega)}(\mathsf{C}))}\).

Claim~\ref{CGdwUfwE} follows.

We next prove that the morphisms \(\alpha^{[\bigcdot]}=\left(\alpha^{[\mathsf{C}]}\right)_{\mathsf{C}\in\mathrm{Ob}(\w\mathsf{Cat}^{\mathrm{MS}})}\), when \(\mathsf{C}\) ranges over the objects of \(\w\mathsf{Cat}^{\mathrm{MS}}\), are the components of a natural transformation.

\begin{claim}\label{CGdwUfwF}
For every many-sorted \(\omega\)-functor \(F:\mathsf{C}\mor\mathsf{D}\), the square of many-sorted \(\omega\)-functors in Figure~\ref{FEtaNatTrans2} commutes.
\begin{figure}[t]
\[
\xymatrix{
\mathsf{C}
  \ar[d]_{F}
  \ar[r]^-{\alpha^{[\mathsf{C}]}}
&
\mathrm{Gd}^{(\omega)}(\mathrm{Uf}^{(\omega)}(\mathsf{C}))
  \ar[d]^{\mathrm{Gd}^{(\omega)}(\mathrm{Uf}^{(\omega)}(F))}
\\
\mathsf{D}
  \ar[r]_-{\alpha^{[\mathsf{D}]}}
&
\mathrm{Gd}^{(\omega)}(\mathrm{Uf}^{(\omega)}(\mathsf{D}))
}
\]
\caption{The natural transformation \(\alpha^{[\mathsf{C}]}\) on many-sorted \(\omega\)-functors.}
\label{FEtaNatTrans2}
\end{figure}
\end{claim}

Let \(\mathsf{C}\) and \(\mathsf{D}\) be the \(\omega\)-categories under consideration, where
\begin{align*}
\mathsf{C}&=
\left(
(C_{k})_{k\in\omega},
\left(
\circ^{(j)\mathsf{C}},\mathrm{sc}^{(j,k)\mathsf{C}},\mathrm{tg}^{(j,k)\mathsf{C}},\mathrm{i}^{(k,j)\mathsf{C}}
\right)_{(j,k)\in\coprod_{k\in\omega}k}
\right) 
\mbox{ and}
\\
\mathsf{D}&=
\left(
(D_{k})_{k\in\omega},
\left(
\circ^{(j)\mathsf{D}},\mathrm{sc}^{(j,k)\mathsf{D}},\mathrm{tg}^{(j,k)\mathsf{D}},\mathrm{i}^{(k,j)\mathsf{D}}
\right)_{(j,k)\in\coprod_{k\in\omega}k}
\right)
\end{align*}
and let
\begin{align*}
\mathrm{Gd}^{(\omega)}\left(\mathrm{Uf}^{(\omega)}(\mathsf{C})\right)
&=
\left(
\left(
[C_{\omega}]_{k}\right)_{k\in\omega},
\left(
\circ'^{(j)\mathsf{C}},\mathrm{sc}'^{(j,k)\mathsf{C}},\mathrm{tg}'^{(j,k)\mathsf{C}},\mathrm{i}'^{(k,j)\mathsf{C}}
\right)_{(j,k)\in\coprod_{k\in\omega}k}
\right)
\mbox{ and}
\\
\mathrm{Gd}^{(\omega)}\left(\mathrm{Uf}^{(\omega)}(\mathsf{D})\right)
&=
\left(
\left(
[D_{\omega}]_{k}\right)_{k\in\omega},
\left(
\circ'^{(j)\mathsf{D}},\mathrm{sc}'^{(j,k)\mathsf{D}},\mathrm{tg}'^{(j,k)\mathsf{D}},\mathrm{i}'^{(k,j)\mathsf{D}}
\right)_{(j,k)\in\coprod_{k\in\omega}k}
\right)
\end{align*}
be the corresponding graduations of the unifications of \(\mathsf{C}\) and \(\mathsf{D}\), respectively.

Let \(F=(F_{k})_{k\in\omega}\) be a many-sorted \(\omega\)-functor from \(\mathsf{C}\) to \(\mathsf{D}\). Let us recall that according to Definition~\ref{DUfwFun}, its unification is the single-sorted \(\omega\)-functor that assigns to every \(\omega\)-cell class \([f]_{j}\) in \([C_{\omega}]\) the \(\omega\)-cell class \([F_{j}(f)]_{j}\) in \([D_{\omega}]\).  Moreover, according to Remark~\ref{RGdw}, the graduation of the unification of \(F\) is given by
\[
\mathrm{Gd}^{(\omega)}\left(\mathrm{Uf}^{(\omega)}(F)\right)
=
\left(
\mathrm{Uf}^{(\omega)}(F)\big|_{[C_{\omega}]_{k}}^{[D_{\omega}]_{k}}
\right)_{k\in\omega}.
\]

Let \(\alpha^{[\mathsf{C}]}\) be the many-sorted \(\omega\)-functor from \(\mathsf{C}\) to \(\mathrm{Gd}^{(\omega)}(\mathrm{Uf}^{(\omega)}(\mathsf{C}))\) and \(\alpha^{[\mathsf{D}]}\) be the many-sorted \(\omega\)-functor from \(\mathsf{D}\) to \(\mathrm{Gd}^{(\omega)}(\mathrm{Uf}^{(\omega)}(\mathsf{D}))\) defined before. Let us recall that, for every \(i\in\omega\) and every \(f\in C_{i}\) and \(g\in C_{i}\), \(\alpha^{[\mathsf{C}]}_{i}\) and \(\alpha^{[\mathsf{D}]}_{i}\) are defined as
\begin{align}
\alpha^{[\mathsf{C}]}_{i}(f)&=[f]_{i}
&\alpha^{[\mathsf{D}]}_{i}(g)&=[g]_{i}
\end{align}

We represent in the  diagram of Figure~\ref{FEtaNTwCat} all the different sets and mappings under consideration.

\begin{figure}
\begin{center}
\[
\xymatrixcolsep={4.5ex}
\xymatrixrowsep={4.5ex}
\xymatrix@!0{
  &&&&&&&&
  \dots
    \ar@<1ex>[ddll]
    \ar@<-1ex>[ddll]
  &&&&&&&&&
  \dots
    \ar@<1ex>[ddll]
    \ar@<-1ex>[ddll]
\\
  &&&&&&&&&&&&&&&&&
\\
  &&&&&&
  C_{n}
    \ar[rrrrrrrrr]|(.35){\scriptstyle \alpha^{[\mathsf{C}]}_{n}}
    \ar[ddddd]|(.6){\scriptstyle F_{n}}|(.8)\hole
    \ar[rruu]
    \ar@<1ex>[ddll]
    \ar@<-1ex>[ddll]
  &&&&&&&&&
  [C_{\omega}]_{n}
    \ar[ddddd]|{\scriptstyle\mathrm{Gd}^{(\omega)}(\mathrm{Uf}^{(\omega)}(F))_{n}}
    \ar[rruu]
    \ar@<1ex>[ddll]
    \ar@<-1ex>[ddll]
  &&
\\
  &&&&&&&&&&&&&&&&&
\\
  &&&&
  \dots
    \ar[rruu]
    \ar@<1ex>[ddll]
    \ar@<-1ex>[ddll]
  &&&&&&&&&
  \dots
    \ar[rruu]
    \ar@<1ex>[ddll]
    \ar@<-1ex>[ddll]
  &&&&
\\
  &&&&&&&&
  \dots
    \ar@<1ex>[ddll]|(.43)\hole
    \ar@<-1ex>[ddll]|(.57)\hole
  &&&&&&&&&
  \dots
    \ar@<1ex>[ddll]
    \ar@<-1ex>[ddll]
\\
  &&
  C_{1}
    \ar[rrrrrrrrr]|(.35){\scriptstyle \alpha^{[\mathsf{C}]}_{1}}
    \ar[ddddd]|(.4)\hole|(.6){\scriptstyle F_{1}}
    \ar[rruu]
    \ar@<1ex>[ddll]
    \ar@<-1ex>[ddll]
  &&&&&&&&&
  [C_{\omega}]_{1}
    \ar[ddddd]|{\scriptstyle\mathrm{Gd}^{(\omega)}(\mathrm{Uf}^{(\omega)}(F))_{1}}
    \ar[rruu]
    \ar@<1ex>[ddll]
    \ar@<-1ex>[ddll]
  &&&&&&
\\
  &&&&&&
  D_{n}
    \ar[rrrrrrrrr]|(.35){\scriptstyle \alpha^{[\mathsf{D}]}_{n}}|(0.42)\hole|(0.45)\hole|(0.48)\hole|(.555)\hole
    \ar[rruu]|(0.5)\hole
    \ar@<1ex>[ddll]|(0.43)\hole
    \ar@<-1ex>[ddll]|(0.58)\hole
  &&&&&&&&&
  [D_{\omega}]_{n}
    \ar[rruu]
    \ar@<1ex>[ddll]
    \ar@<-1ex>[ddll]
  &&
\\
  C_{0}
    \ar[rrrrrrrrr]|(.35){\scriptstyle \alpha^{[\mathsf{C}]}_{0}}
    \ar[ddddd]|(.6){\scriptstyle F_{0}}
    \ar[rruu]
  &&&&&&&&&
  [C_{\omega}]_{0}
    \ar[ddddd]|{\scriptstyle\mathrm{Gd}^{(\omega)}(\mathrm{Uf}^{(\omega)}(F))_{0}}
    \ar[rruu]
  &&&&&&&&
\\
  &&&&
  \dots
    \ar[rruu]|(0.5)\hole
    \ar@<1ex>[ddll]
    \ar@<-1ex>[ddll]
  &&&&&&&&&
  \dots
    \ar[rruu]
    \ar@<1ex>[ddll]
    \ar@<-1ex>[ddll]
  &&&&
\\
  &&&&&&&&&&&&&&&&&
\\
  &&
  D_{1}
    \ar[rrrrrrrrr]|(.35){\scriptstyle \alpha^{[\mathsf{D}]}_{1}}|(.78)\hole
    \ar[rruu]
    \ar@<1ex>[ddll]
    \ar@<-1ex>[ddll]
  &&&&&&&&&
  [D_{\omega}]_{1}
    \ar[rruu]
    \ar@<1ex>[ddll]
    \ar@<-1ex>[ddll]
  &&&&&&
\\
  &&&&&&&&&&&&&&&&&
\\
  D_{0}
    \ar[rrrrrrrrr]|(.35){\scriptstyle \alpha^{[\mathsf{D}]}_{0}}
    \ar[rruu]
  &&&&&&&&&
  [D_{\omega}]_{0}
    \ar[rruu]
  &&&&&&&& 
}
\]
\end{center}
\caption{The natural transformation $\alpha^{[\bigcdot]}$ on many-sorted $\omega$-categories.}
\label{FEtaNTwCat}
\end{figure}

In order to prove Claim~\ref{CGdwUfwF} we need to check that, for every $i\in\omega$, the diagram in Figure~\ref{FEtaNTi} commutes.

\begin{figure}
\begin{center}
\begin{tikzpicture}
[ACliment/.style={-{To [angle'=45, length=5.75pt, width=4pt, round]}}
, scale=1, 
AClimentD/.style={dogble eqgal sign distance,
-implies
}
]

\node[] (a) at (0,0) [] {$C_{i}$};
\node[] (b) at (3,0) [] {$\mathrm{Gd}^{(\omega)}(\mathrm{Uf}^{(\omega)}(\mathsf{C}))_{i}$};
\node[] (c) at (0,-1.5) [] {$D_{i}$};
\node[] (d) at (3,-1.5) [] {$\mathrm{Gd}^{(\omega)}(\mathrm{Uf}^{(\omega)}(\mathsf{D}))_{i}$};

\draw[ACliment] (a) to node [above] {$\scriptstyle \alpha^{[\mathsf{C}]}_{i}$} (b);
\draw[ACliment] (a) to node [left] {$\scriptstyle F_{i}$} (c);
\draw[ACliment] (c) to node [below] {$\scriptstyle \alpha^{[\mathsf{D}]}_{i}$} (d);
\draw[ACliment] (b) to node [right] {$\scriptstyle \mathrm{Gd}^{(\omega)}(\mathrm{Uf}^{(\omega)}(F))_{i}$} (d);
\end{tikzpicture}
\end{center}
\caption{The natural transformation $\eta$ at layer $i$ of a many-sorted $\omega$-category.}
\label{FEtaNTi}
\end{figure}

Let \(i\in\omega\) and \(f\in C_{i}\), the the following chain of equalities holds
\allowdisplaybreaks
\begin{align*}
\left(\mathrm{Gd}^{(\omega)}\circ\mathrm{Uf}^{(\omega)}(F)\right)_{i}(\alpha^{[\mathsf{C}]}_{i}(f))
&=
\left(\mathrm{Gd}^{(\omega)}\circ\mathrm{Uf}^{(\omega)}(F)\right)_{i}\left(\left[f\right]_{i}\right)
\tag{1}
\\
&=
\left(\mathrm{Uf}^{(\omega)}(F)\big|^{[D_{\omega}]_{i}}_{[C_{\omega}]_{i}}\right)\left(\left[f\right]_{i}\right)
\tag{2}
\\
&=
\mathrm{Uf}^{(\omega)}(F)\left(\left[f\right]_{i}\right)
\tag{3}
\\
&=
\left[F_{i}(f)\right]_{i}
\tag{4}
\\
&=
\alpha^{[\mathsf{D}]}_{i}\left(F_{i}(f)\right).
\tag{5}
\end{align*}

The first equality unravels the \(i\)-th component of the many-sorted \(\omega\)-functor \(\alpha^{[\mathsf{C}]}\); the second equality unravels the \(i\)-th component of the many-sorted \(\omega\)-functor \(\mathrm{Gd}^{(\omega)}\circ\mathrm{Uf}^{(\omega)}(F)\); the third equality applies the definition of restriction of a mapping; the fourth equality applies the definition of the functor \(\mathrm{Uf}^{(\omega)}(F)\); finally, the last equality recovers the \(i\)-th component of the many-sorted \(\omega\)-functor \(\alpha^{[\mathsf{D}]}\)

Claim~\ref{CGdwUfwF} follows.

So, considering the foregoing, we can affirm that \(\alpha^{[\bigcdot]}\) is a natural transformation from \(\mathrm{Id}^{\w\mathsf{Cat}^{\mathrm{MS}}}\) to \(\mathrm{Gd}^{(\omega)}\circ\mathrm{Uf}^{(\omega)}\).

Let us note that each many-sorted \(\omega\)-category \(\mathsf{C}\) in \(\w\mathsf{Cat}^{\mathrm{MS}}\) determines, by Claim~\ref{CGdwUfwC}, the many-sorted \(\omega\)-functor \(\alpha^{[\mathsf{C}]}\) which, by Claim~\ref{CGdwUfwE}, is an isomorphism.

Thus, \(\alpha^{[\bigcdot]}:\mathrm{Id}^{\w\mathsf{Cat}^{\mathrm{MS}}}\cel\mathrm{Gd}^{(\omega)}\circ\mathrm{Uf}^{(\omega)}\) is a natural isomorphism.
\end{proof}

From the above, Proposition~\ref{Eqvt} and Definition~\ref{lali} we obtain the following corollary.

\begin{corollary}\label{CDwCatEquiv} 
The functor $\mathrm{Uf}^{(\omega)}$ is a left adjoint left inverse of the functor $\mathrm{Gd}^{(\omega)}$. Hence the categories $\w\mathsf{Cat}$ and $\w\mathsf{Cat}^{\mathrm{MS}}$ are equivalent.
\end{corollary}

\section{
\texorpdfstring
{$n$-categorial $\Sigma$-algebras}
{n-categorial algebras}
}

In this paper, occasionally, we will use a slight generalization of the notion of many-sorted algebra, which we define as follows. 

\begin{definition}\label{DnCat}
We will call the objects of the functor category $\mathsf{nCat}^{S}$ \emph{$S$-sorted $n$-categories}. Thus, an $S$-sorted $n$-category is an $S$-indexed family $\mathsf{A}=(\mathsf{A}_{s})_{s\in S}$ such that, for every sort $s\in S$, $\mathsf{A}_{s}$ is an $n$-category. Given two $S$-sorted $n$-categories 
$\mathsf{A}$ and $\mathsf{B}$, we will call the morphisms in $\mathsf{nCat}^{S}$ from $\mathsf{A}$ to 
$\mathsf{B}$ \emph{$S$-sorted $n$-functors}. Thus, an $S$-sorted $n$-functor from $\mathsf{A}$ to $\mathsf{B}$ is an $S$-indexed family $F=(F_{s})_{s\in S}$ such that, for every sort $s\in S$, $F_{s}$ is an $n$-functor from $\mathsf{A}_{s}$ to $\mathsf{B}_{s}$. We will denote by $\mathrm{nFunc}(\mathsf{A},\mathsf{B})$ the set of all $S$-sorted $n$-functors from $\mathsf{A}$ to $\mathsf{B}$. We consider the expressions $F\in\mathrm{nFunc}(\mathsf{A},\mathsf{B})$ and $F\colon\mathsf{A}\mor\mathsf{B}$ as synonyms.
We recall that for $F\colon \mathsf{A}\mor\mathsf{B}$ and $G\colon \mathsf{B}\mor\mathsf{C}$, the \emph{composition} $G\circ F$ is $(G_{s}\circ F_{s})_{s\in S}$, and that for an $S$-sorted $n$-category $\mathsf{A}$, the \emph{identity} 
$S$-sorted $n$-functor at $\mathsf{A}$ $\mathrm{Id}_{\mathsf{A}}$ is $(\mathrm{Id}_{\mathsf{A}_{s}})_{s\in S}$. 
\end{definition}

%

\begin{definition}\label{DnCatAlg}
The $S^{\star}\times S$-sorted set of the finitary operations on an $S$-sorted $n$-category $\mathsf{A}$ is $(\mathrm{nFunc}(\mathsf{A}_{\mathbf{s}},\mathsf{A}_{s}))_{(\mathbf{s},s)\in S^{\star}\times S}$, where, for every $\mathbf{s}\in S^{\star}$, $\mathsf{A}_{\mathbf{s}}=\prod_{j\in\bb{\mathbf{s}}}\mathsf{A}_{s_{j}}$, with $\bb{\mathbf{s}}$ denoting the length of the word $\mathbf{s}$ (if $\mathbf{s}=\lambda$, then $\mathsf{A}_{\lambda}$ is a final $n$-category). 

Let $\Sigma$ be an $S$-sorted signature. A \emph{structure of $n$-categorial $\Sigma$-algebra on} an $S$-sorted $n$-category $\mathsf{A}$ is a family $F = (F_{\mathbf{s},s})_{(\mathbf{s},s)\in S^{\star}\times S}$, where, for $(\mathbf{s},s)\in S^{\star}\times S$, $F_{\mathbf{s},s}$ is a mapping from $\Sigma_{\mathbf{s},s}$ to $\mathrm{nFunc}(\mathsf{A}_{\mathbf{s}},\mathsf{A}_{s})$ (if $(\mathbf{s},s) = (\lambda,s)$ and $\sigma\in \Sigma_{\lambda,s}$, then $F_{\lambda,s}(\sigma)$, picks out an object of $\mathsf{A}_{s}$ and its identity morphism).  An $S$-\emph{sorted $n$-categorial $\Sigma$-algebra} is a pair $(\mathsf{A},F)$, abbreviated to $\boldsymbol{\mathsf{A}}$, where $\mathsf{A}=(\mathsf{A}_{s})_{s\in S}$ is an $S$-sorted $n$-category and $F$ a structure of $n$-categorial $\Sigma$-algebra on $\mathsf{A}$.  For a pair $(\mathbf{s},s)\in S^{\star}\times S$ and a formal operation $\sigma\in \Sigma_{\mathbf{s},s}$, in order to simplify the notation, the $n$-functor $F_{\mathbf{s},s}(\sigma)$ from $\mathsf{A}_{\mathbf{s}}$ to $\mathsf{A}_{s}$ will be written as $F_{\sigma}$, or simply as $\sigma^{\boldsymbol{\mathsf{A}}}$. 

An $S$-\emph{sorted $n$-categorial $\Sigma$-homomorphism} from $\boldsymbol{\mathsf{A}}$ to $\boldsymbol{\mathsf{B}}$, where $\boldsymbol{\mathsf{B}}=(\mathsf{B},G)$, is a triple $(\boldsymbol{\mathsf{A}},F,\boldsymbol{\mathsf{B}})$, abbreviated to $F\colon\boldsymbol{\mathsf{A}}\mor\boldsymbol{\mathsf{B}}$, where $F$ is an $S$-sorted $n$-functor from $\mathsf{A}$ to $\mathsf{B}$ such that, for every $k\in n$, every $(\mathbf{s},s)\in S^{\star}\times S$, every $\sigma\in\Sigma_{\mathbf{s},s}$ and every family of $k$-cells $(\mathsf{a}_{j})_{j\in\bb{\mathbf{s}}}\in\mathsf{A}_{\mathbf{s}}$, we have that
$$
F_{s}\left(\sigma^{\boldsymbol{\mathsf{A}}}
\left(\left(\mathsf{a}_{j}
\right)_{j\in\bb{\mathbf{s}}}\right)\right)
=
\sigma^{\boldsymbol{\mathsf{B}}}
\left(\left(
F_{s_{j}}\left(
\mathsf{a}_{j}
\right)
\right)_{j\in\bb{\mathbf{s}}}
\right).
$$
We denote by $\mathsf{nCat}^{S}\mathsf{Alg}(\Sigma)$ the category of $S$-sorted $n$-categorial  $\Sigma$-algebras and $S$-sorted $n$-categorial $\Sigma$-homomorphisms.
\end{definition}

\begin{remark}
The above definitions can be extended from $\mathsf{nCat}^{S}$, the category of $S$-sorted $n$-categories and $S$-sorted $n$-functors, up to the 2-category of $S$-sorted $n$-categories, $S$-sorted $n$-functors and $S$-sorted $n$-natural transformations. We leave this possibility open for future work.
\end{remark}


\part{First-order rewriting systems}
\chapter{
\texorpdfstring
{Paths on terms}
{Paths}
}\label{S1A}

In this chapter we begin by defining the notion of many-sorted term rewriting system. After that, for a many-sorted term rewriting system $\boldsymbol{\mathcal{A}}$, we define the concept of path in $\boldsymbol{\mathcal{A}}$ from a term to another. Then we consider the many-sorted set of paths in $\boldsymbol{\mathcal{A}}$, denoted by $\mathrm{Pth}_{\boldsymbol{\mathcal{A}}}$. We further define the notions of $(0,1)$-source of a path, $(0,1)$-target of a path and $(1,0)$-identity path, denoted by $\mathrm{sc}^{(0,1)}$, $\mathrm{tg}^{(0,1)}$ and $\mathrm{ip}^{(1,0)\sharp}$, respectively. We then start the structural study of these objects by characterising the $(1,0)$-identity paths. We also define the partial operation of $0$-composition of paths and prove that it is well-defined. Following this we define the notion of subpath of a path and we relate the concepts of $0$-composition and subpaths. We then associate to every rewrite rule a special path which we call an echelon and we define the concept of echelonless path. Next we characterise the echelons by means of the translations that occur along a path. We then define the notion of head-constant echelonless path and prove that echelonless paths are head-constant. We conclude this chapter by introducing an algorithm which we call the path extraction algorithm for echelonless paths.


\begin{definition}
\label{DRewSys}
A \emph{many-sorted term rewriting system}\index{rewriting system!first-order, $\boldsymbol{\mathcal{A}}$} (or, simply, a \emph{rewriting system}) is an ordered quadruple
$
(S,\Sigma,X,\mathcal{A}),
$
often abbreviated to $\boldsymbol{\mathcal{A}}$, 
where $S$ is a set of sorts, $\Sigma$ an $S$-sorted signature, $X$ an $S$-sorted set and, for the many-sorted signature $\mathbf{\Sigma} = (S,\Sigma)$, $\mathcal{A}$ a subset of $\mathrm{Rwr}(\mathbf{\Sigma}, X) = (\mathrm{T}_{\Sigma}(X)^{2}_{s})_{s\in S}$, the $S$-sorted set of the \emph{rewrite rules with variables in $X$}, where $\mathrm{T}_{\Sigma}(X)$ is the underlying $S$-sorted set of $\mathbf{T}_{\Sigma}(X)$, the free many-sorted $\Sigma$-algebra on $X$. 

\index{rewrite rule!first-order, $\mathfrak{p}$}
For $s\in S$,  we will call the elements of $\mathrm{Rwr}(\mathbf{\Sigma}, X)_{s}$ \emph{rewrite rules of type $(X,s)$} and we will denote them with lowercase Euler fraktur letters, with or without subscripts, e.g., $\mathfrak{p}$, $\mathfrak{p}_{i}$, $\mathfrak{q}$, $\mathfrak{q}_{i}$, etc). 

We will say that $\boldsymbol{\mathcal{A}}$ is a \emph{finite} rewriting system if $\mathcal{A}$ is finite, i.e., if $\mathrm{supp}_{S}(\mathcal{A})$ is finite and, for every $s\in \mathrm{supp}_{S}(\mathcal{A})$, $\mathcal{A}_{s}$ is finite.
\end{definition}

\begin{remark}
For a many-sorted signature $\mathbf{\Sigma} = (S,\Sigma)$, an $S$-sorted set $X$ and a sort $s\in S$, a 
$\mathbf{\Sigma}$-equation of type $(X,s)$, or a $\mathbf{\Sigma}$-equation of sort $s$ with variables in $X$, is an ordered pair $(P,Q)\in \mathrm{T}_{\Sigma}(X)^{2}_{s}$ and the $S$-sorted set of the equations with variables in $X$, denoted by $\mathrm{Eqt}(\mathbf{\Sigma})_{X}$,  is $(\mathrm{T}_{\Sigma}(X)^{2}_{s})_{s\in S}$, which is $\mathrm{Rwr}(\mathbf{\Sigma},X)$. So, what is the difference between equations and rewrite rules? Well, the difference between the two notions is given by the following considerations.     
Formally, one obtains a rewrite rule with variables in $X$ from an equation with variables in $X$ by simply orienting the equation (or, swapping roles, one obtains an equation by seeing it as a two-way rewrite rule). This will be notationally expressed by writing equations as $P = Q$ and rewrite rules as $(P,Q)$. By doing so, a rewrite rule $(P,Q)$ is considered to be oriented in the sense that an occurrence of $P$ (in a term) may be rewritten as $Q$ but not vice versa, and it gives rise to a translation operator. On the other hand, an equation $P = Q$ serves to determine the set of all those algebras in which the terms of the equation, realized as operations on them, are identical, i.e., the set $\mathrm{Mod}(\{P = Q\}) = \{\mathbf{A}\in\mathrm{Alg}(\Sigma)\mid \mathbf{A}\models P = Q\}$ of all models of $P = Q$.
\end{remark}

\begin{assumption}
From now on $\boldsymbol{\mathcal{A}}$ stands for a rewriting system, fixed once and for all.
\end{assumption}

We next define the notion of path in $\boldsymbol{\mathcal{A}}$ from a term to another.

\begin{restatable}{definition}{DPth}
\label{DPth}
Let $s$ be a sort in $S$, $P$, $Q$ terms in $\mathrm{T}_{\Sigma}(X)_{s}$, and $\mathbf{c}$ a word in $S^{\star}$. Then a $\mathbf{c}$-\emph{path in $\boldsymbol{\mathcal{A}}$ from $P$ to $Q$} is an ordered triple
\index{path!first-order!$\mathfrak{P}$}
\begin{equation*}
\mathfrak{P} =
\left(
(P_{i})_{i\in \bb{\mathbf{c}}+1},
(\mathfrak{p}_{i})_{i\in\bb{\mathbf{c}}},
(T_{i})_{i\in\bb{\mathbf{c}}}
\right)
\end{equation*}
in
$
\mathrm{T}_{\Sigma}(X)_{s}^{\bb{\mathbf{c}}+1}
\times
\mathcal{A}_{\mathbf{c}}
\times 
\mathrm{Tl}_{\mathbf{c}}(\mathbf{T}_{\Sigma}(X))_{s},
$ 
where $\mathcal{A}_{\mathbf{c}} = \prod_{i\in\bb{\mathbf{c}}}\mathcal{A}_{c_{i}}$ and 
$\mathrm{Tl}_{\mathbf{c}}(\mathbf{T}_{\Sigma}(X))_{s} =
\prod_{i\in\bb{\mathbf{c}}}\mathrm{Tl}_{c_{i}}(\mathbf{T}_{\Sigma}(X))_{s}$, with $\mathrm{Tl}_{c_{i}}(\mathbf{T}_{\Sigma}(X))_{s}$ the set of the $c_{i}$-translations of sort $s$ for $\mathbf{T}_{\Sigma}(X)$, see Definition~\ref{DTrans}, such that

\begin{enumerate}
\item $P_{0}=P$,
\item $P_{\bb{\mathbf{c}}}=Q$, and,
\item for every $i\in\bb{\mathbf{c}}$, if $\mathfrak{p}_{i}=(M_{i},N_{i})$, then
\begin{enumerate}
\item[(i)] 
$T_{i}(M_{i})=P_{i}$ \text{ and}
\item[(ii)] 
$T_{i}(N_{i})=P_{i+1}$.
\end{enumerate}
\end{enumerate}
That is, at each step $i\in\bb{\mathbf{c}}$, we consider a rewrite rule of type $(X,c_{i})$, $\mathfrak{p}_{i}$, and a $c_{i}$-translation of sort $s$ for $\mathbf{T}_{\Sigma}(X)$, $T_{i}$, and we require that (i) the translation by $T_{i}$ of $M_{i}$ is $P_{i}$, whilst (ii) the translation by $T_{i}$ of $N_{i}$ is $P_{i+1}$. This statement can also be understood through the use of substitutions. We will be say that $P_{i}$ contains $M_{i}$ as a subterm and that $P_{i+1}$ results from substituting one of its subterms $M_{i}$ for $N_{i}$ in $P_{i}$. This justifies the name rewriting rules for the elements in the family 
$(\mathfrak{p}_{i})_{i\in\bb{\mathbf{c}}}$. On the other hand, we could think of the translations in the family $(T_{i})_{i\in\bb{\mathbf{c}}}$ as the contexts in which the rewriting rules are applied.

These paths will be variously depicted as $\mathfrak{P}\colon P\mor Q$, $\mathfrak{P}\colon P_{0}\mor P_{\bb{\mathbf{c}}}$,  or
$$
\xymatrix@C=55pt{
\mathfrak{P}: P_{0}
\ar[r]^-{\text{\Small{($\mathfrak{p}_{0}$, $T_{0}$)}}}
&
P_{1}
\ar[r]^-{\text{\Small{($\mathfrak{p}_{1}$, $T_{1}$)}}}
&
{}
\hdots
{}
\ar[r]^-{\text{\Small{($\mathfrak{p}_{\bb{\mathbf{c}}-2}$, $T_{\bb{\mathbf{c}}-2}$, )}}}
&
P_{\bb{\mathbf{c}}-1}
\ar[r]^-{\text{\Small{($\mathfrak{p}_{\bb{\mathbf{c}}-1}$, $T_{\bb{\mathbf{c}}-1}$,)}}}
&
P_{\bb{\mathbf{c}}}
}
$$
(the following variant
$$
\xymatrix@C=55pt{
\mathfrak{P}: P_{0}
\ar[r]^-{\text{\Small{$(\mathfrak{p}$, $T)_{0}$}}}
&
P_{1}
\ar[r]^-{\text{\Small{$(\mathfrak{p}$, $T)_{1}$}}}
&
{}
\hdots
{}
\ar[r]^-{\text{\Small{$(\mathfrak{p}$, $T)_{\bb{\mathbf{c}}-2}$}}}
&
P_{\bb{\mathbf{c}}-1}
\ar[r]^-{\text{\Small{$(\mathfrak{p}$, $T)_{\bb{\mathbf{c}}-1}$}}}
&
P_{\bb{\mathbf{c}}}
}
$$
will be occasionally used).

For every $i\in \bb{\mathbf{c}}$, we will say that $P_{i+1}$ is
$(\mathfrak{p}_{i}, T_{i})$-\emph{directly derivable}\index{directly derivable} or, when no confusion can arise, \emph{directly derivable} from $P_{i}$. For every $i\in \bb{\mathbf{c}}+1$, the term $P_{i}$ will be called a \emph{$0$-constituent}\index{constituent!$0$-constituent} of the $\bb{\mathbf{c}}$-path $\mathfrak{P}$. Let us note that all $0$-constituents of a path $\mathfrak{P}$ have the same sort. The term $P_{0}$ will be called the $(0,1)$-\emph{source}\index{source!first-order!$\mathrm{sc}^{(0,1)}$} of the $\mathbf{c}$-path $\mathfrak{P}$, the term $P_{\bb{\mathbf{c}}}$ will be called the $(0,1)$-\emph{target}\index{target!first-order!$\mathrm{tg}^{(0,1)}$} of the $\mathbf{c}$-path $\mathfrak{P}$, and we will say that $\mathfrak{P}$ is a \emph{path from} $P_{0}$ \emph{to} $P_{\bb{\mathbf{c}}}$.

The \emph{length}\index{path!first-order!length} of a $\mathbf{c}$-path $\mathfrak{P}$ in $\boldsymbol{\mathcal{A}}$, denoted by $\bb{\mathfrak{P}}$, is $\bb{\mathbf{c}}$ and we will say that $\mathfrak{P}$ has $\bb{\mathbf{c}}$ \emph{steps}\index{path!first-order!step}. 
If $\bb{\mathfrak{P}}=0$, then we will say that $\mathfrak{P}$ is a  \emph{$(1,0)$-identity path}\index{identity!first-order!$\mathrm{ip}^{(0,1)\sharp}$} if, and only if, there exists a sort $s$ in $S$ and a term $P$ in $\mathrm{T}_{\Sigma}(X)_{s}$ such that, for $\lambda\in S^{\star}$, the empty word on $S$, we have that $((P),\lambda,\lambda)$, identified to $(P,\lambda,\lambda)$, where, by abuse of notation, we have written $(\lambda,\lambda)$ for the unique element of $\mathcal{A}_{\lambda}\times \mathrm{Tl}_{\lambda}(\mathbf{T}_{\Sigma}(X))_{s}$, is equal to $\mathfrak{P}$. This path will be called the \emph{$(1,0)$-identity path on $P$}. 

If $\bb{\mathfrak{P}}=1$, then we will say that $\mathfrak{P}$ is a \emph{one-step}\index{path!first-order!one-step} path. Moreover, in this case, i.e., when $\mathbf{c} = (c)$, for a unique $c\in S$, and identifying, when no confusion can arise, $(c)$ with $c$, we will speak of $c$-paths, instead of $(c)$-paths.

We will denote by
\begin{enumerate}
\item $\mathrm{Pth}_{\mathbf{c},\boldsymbol{\mathcal{A}},s}(P,Q)$ the set of all $\mathbf{c}$-paths in $\boldsymbol{\mathcal{A}}$ from $P$ to $Q$; by
\item $\mathrm{Pth}_{\boldsymbol{\mathcal{A}},s}(P,Q)$ the set $\bigcup_{\mathbf{c}\in S^{\star}}\mathrm{Pth}_{\mathbf{c},\boldsymbol{\mathcal{A}},s}(P,Q)$ and we will call its elements \emph{paths in $\boldsymbol{\mathcal{A}}$ from $P$ to $Q$};
by
\item $\mathrm{Pth}_{\boldsymbol{\mathcal{A}},s}(P,\cdot)$ the set $\bigcup_{Q\in \mathrm{T}_{\Sigma}(X)_{s}}\mathrm{Pth}_{\boldsymbol{\mathcal{A}},s}(P,Q)$ and we will call its elements \emph{paths in $\boldsymbol{\mathcal{A}}$ from $P$};
by
\item $\mathrm{Pth}_{\boldsymbol{\mathcal{A}},s}(\cdot,Q)$ the set $\bigcup_{P\in \mathrm{T}_{\Sigma}(X)_{s}}\mathrm{Pth}_{\boldsymbol{\mathcal{A}},s}(P,Q)$ and we will call its elements \emph{paths in $\boldsymbol{\mathcal{A}}$ to $Q$}; by
\item $\mathrm{Pth}_{\boldsymbol{\mathcal{A}},s}$ the set $\bigcup_{P,Q\in \mathrm{T}_{\Sigma}(X)_{s}}\mathrm{Pth}_{\boldsymbol{\mathcal{A}},s}(P,Q)$; and, finally, by
\item $\mathrm{Pth}_{\boldsymbol{\mathcal{A}}}$ the $S$-sorted set $(\mathrm{Pth}_{\boldsymbol{\mathcal{A}},s})_{s\in S}$.
\end{enumerate}

We will denote by
\begin{enumerate}
\item $\mathrm{ip}^{(1,X)}$\index{identity!first-order!$\mathrm{ip}^{(1,X)}$} the $S$-sorted mapping from $X$ to $\mathrm{Pth}_{\boldsymbol{\mathcal{A}}}$ that, for every sort $s\in S$, sends $x\in X_{s}$ to $(x,\lambda,\lambda)$, the $(1,0)$-identity path on $x$; by
\item $\mathrm{sc}^{(0,1)}$ the $S$-sorted mapping from $\mathrm{Pth}_{\boldsymbol{\mathcal{A}}}$ to $\mathrm{T}_{\Sigma}(X)$ that sends a path to its $(0,1)$-source; by
\item $\mathrm{tg}^{(0,1)}$ the $S$-sorted mapping from $\mathrm{Pth}_{\boldsymbol{\mathcal{A}}}$ to $\mathrm{T}_{\Sigma}(X)$ that sends a path to its $(0,1)$-target; and by
\item $\mathrm{ip}^{(1,0)\sharp}$ the $S$-sorted mapping that sends a term $P$ to $(P,\lambda,\lambda)$, the $(1,0)$-identity path on $P$. 
\end{enumerate}
These $S$-sorted mappings are depicted in the diagram of Figure~\ref{FPthX}.

\index{path!first-order!$\mathrm{Pth}_{\boldsymbol{\mathcal{A}}}$}
 
Finally, given a sort $s$ in $S$ and a path $\mathfrak{P}$ in $\mathrm{Pth}_{\boldsymbol{\mathcal{A}},s}$ we will say that $\mathfrak{P}$ is a \emph{$(0,1)$-loop}\index{loop!first-order!$(0,1)$-loop} if $\mathrm{sc}^{(0,1)}_{s}(\mathfrak{P})=\mathrm{tg}^{(0,1)}_{s}(\mathfrak{P})$. Let us note that every $(1,0)$-identity path is a $(0,1)$-loop.
\end{restatable}

\begin{figure}
\begin{tikzpicture}
[ACliment/.style={-{To [angle'=45, length=5.75pt, width=4pt, round]}},scale=1.1]
\node[] (xoq) at (0,0) [] {$X$};
\node[] (txoq) at (6,0) [] {$\mathrm{T}_{\Sigma}(X)$};
\node[] (p) at (6,-3) [] {$\mathrm{Pth}_{\boldsymbol{\mathcal{A}}}$};
\draw[ACliment]  (xoq) to node [above]
{$\eta^{(0,X)}$} (txoq);
\draw[ACliment, bend right=10]  (xoq) to node [below left] {$\mathrm{ip}^{(1,X)}$} (p);

\node[] (B0) at (6,-1.5)  [] {};
\draw[ACliment]  ($(B0)+(0,1.2)$) to node [above, fill=white] {
$\textstyle \mathrm{ip}^{(1,0)\sharp}$
} ($(B0)+(0,-1.2)$);
\draw[ACliment, bend right]  ($(B0)+(.3,-1.2)$) to node [ below, fill=white] {
$\textstyle \mathrm{tg}^{(0,1)}$
} ($(B0)+(.3,1.2)$);
\draw[ACliment, bend left]  ($(B0)+(-.3,-1.2)$) to node [below, fill=white] {
$\textstyle \mathrm{sc}^{(0,1)}$
} ($(B0)+(-.3,1.2)$);
\end{tikzpicture}
\caption{Many-sorted mappings relative to $X$ at layers 0 \& 1.}\label{FPthX}
\end{figure}

It follows from Definition~\ref{DPth} that the $(1,0)$-identity path on a term is the only path of length $0$ whose $(0,1)$-source (or $(0,1)$-target) equals this term. Consequently, the only interesting properties about $(1,0)$-identity paths are those having to do with their $(0,1)$-source (or $(0,1)$-target). 

We next characterize the $(1,0)$-identity paths as fixed points.

\begin{restatable}{proposition}{PPthId}
\label{PPthId}
Let $s$ be a sort in $S$ and $\mathfrak{P}$ a path in $\mathrm{Pth}_{\boldsymbol{\mathcal{A}},s}$. Then the following statements are equivalent:
\begin{enumerate}
\item $\mathfrak{P}$ is a $(1,0)$-identity path.
\item $\mathfrak{P}=\mathrm{ip}^{(1,0)\sharp}_{s}(\mathrm{sc}^{(0,1)}_{s}(\mathfrak{P}))$, i.e., 
$\mathfrak{P}\in \mathrm{Fix}(\mathrm{ip}^{(1,0)\sharp}_{s}\circ\mathrm{sc}^{(0,1)}_{s})$, the set of fixed points of $\mathrm{ip}^{(1,0)\sharp}_{s}\circ\mathrm{sc}^{(0,1)}_{s}$.
\item $\mathfrak{P}=\mathrm{ip}^{(1,0)\sharp}_{s}(\mathrm{tg}^{(0,1)}_{s}(\mathfrak{P}))$, i.e., 
$\mathfrak{P}\in \mathrm{Fix}(\mathrm{ip}^{(1,0)\sharp}_{s}\circ\mathrm{tg}^{(0,1)}_{s})$, the set of fixed points of $\mathrm{ip}^{(1,0)\sharp}_{s}\circ\mathrm{tg}^{(0,1)}_{s}$.
\end{enumerate}
\end{restatable}

The following are some of the relationships between the $S$-sorted mappings $\eta^{(0,X)}$ (the canonical embedding of $X$ into $\mathrm{T}_{\Sigma}(X)$), $\mathrm{ip}^{(1,X)}$, 
$\mathrm{sc}^{(0,1)}$, $\mathrm{tg}^{(0,1)}$,  $\mathrm{ip}^{(1,0)\sharp}$, and $\mathrm{id}^{\mathrm{T}_{\Sigma}(X)}$ (the identity mapping at $\mathrm{T}_{\Sigma}(X)$). 

\begin{proposition}\label{PBasicEq} The following equalities hold
\begin{enumerate}
\item[(i)] $\mathrm{sc}^{(0,1)}\circ\mathrm{ip}^{(1,0)\sharp}
=\mathrm{id}^{\mathrm{T}_{\Sigma}(X)};$
\item[(ii)] $\mathrm{tg}^{(0,1)}\circ\mathrm{ip}^{(1,0)\sharp}
=\mathrm{id}^{\mathrm{T}_{\Sigma}(X)}.$
\end{enumerate}
The reader is advised to consult the diagram of Figure~\ref{FPthX}.
\end{proposition}

\begin{proposition}\label{PBasicEqX}  The following equalities hold
\begin{itemize}
\item[(i)] $\mathrm{sc}^{(0,1)}\circ\mathrm{ip}^{(1,X)}=\eta^{(0,X)};$
\item[(ii)] $\mathrm{tg}^{(0,1)}\circ\mathrm{ip}^{(1,X)}=\eta^{(0,X)};$
\item[(iii)] $\mathrm{ip}^{(1,0)\sharp}\circ\eta^{(0,X)}=\mathrm{ip}^{(1,X)}.$
\end{itemize}
The reader is advised to consult the diagram of Figure~\ref{FPthX}.
\end{proposition}

\begin{remark}
Later on, after defining a suitable structure of $\Sigma$-algebra on $\mathrm{Pth}_{\boldsymbol{\mathcal{A}}}$, which will give raise to the $\Sigma$-algebra $\mathbf{Pth}^{(0,1)}_{\boldsymbol{\mathcal{A}}}$, we will show that $\mathrm{sc}^{(0,1)}$ and $\mathrm{tg}^{(0,1)}$ are homomorphisms from $\mathbf{Pth}^{(0,1)}_{\boldsymbol{\mathcal{A}}}$ to $\mathbf{T}_{\Sigma}(X)$, and  that $\mathrm{ip}^{(1,0)\sharp}$ is, in fact, the unique homomorphism from $\mathbf{T}_{\Sigma}(X)$ to $\mathbf{Pth}^{(0,1)}_{\boldsymbol{\mathcal{A}}}$ such that $\mathrm{ip}^{(1,0)} = \mathrm{ip}^{(1,0)\sharp}\circ \eta^{(0,X)}$.
\end{remark}

We next define, for every sort $s$ in $S$, the partial operation of $0$-composition of paths of sort $s$.

\begin{restatable}{definition}{DPthComp}
\label{DPthComp}
Let $s$ be a sort in $S$ and $\mathfrak{P}$, $\mathfrak{Q}$ paths in $\mathrm{Pth}_{\boldsymbol{\mathcal{A}},s}$, where, for a unique word $\mathbf{c}\in S^{\star}$, $\mathfrak{P}$ is a path in $\boldsymbol{\mathcal{A}}$ of the form
$$
\mathfrak{P}
=
\left(
(P_{i})_{i\in\bb{\mathbf{c}}+1},
(\mathfrak{p}_{i})_{i\in\bb{\mathbf{c}}},
(T_{i})_{i\in\bb{\mathbf{c}}}
\right),
$$
and, for a unique $\mathbf{d}\in S^{\star}$, $\mathfrak{Q}$ is a $\mathbf{d}$-path in $\boldsymbol{\mathcal{A}}$ of the form
$$
\mathfrak{Q}
=
\left(
(Q_{j})_{j\in\bb{\mathbf{d}}+1},
(\mathfrak{q}_{j})_{j\in\bb{\mathbf{d}}},
(U_{j})_{j\in\bb{\mathbf{d}}}
\right),
$$
such that 
$
\mathrm{sc}^{(0,1)}_{s}(\mathfrak{Q})=\mathrm{tg}^{(0,1)}_{s}(\mathfrak{P}).
$

Then the \emph{$0$-composite}\index{composition!$0$-composition, $\circ^{0}$} of $\mathfrak{P}$ and $\mathfrak{Q}$, denoted by $\mathfrak{Q}\circ^{0}_{s}\mathfrak{P}$, is the ordered triple
$$
\mathfrak{Q}\circ^{0}_{s}\mathfrak{P} =
\left(
(R_{k})_{k\in\bb{\mathbf{e}}+1},
(\mathfrak{r}_{k})_{k\in\bb{\mathbf{e}}},
(V_{k})_{k\in\bb{\mathbf{e}}}
\right),
$$
where $\mathbf{e}=\mathbf{c}\curlywedge\mathbf{d}$, and
\allowdisplaybreaks
\begin{alignat*}{2}
R_{k} &=
\begin{cases}
P_{k}, \\
Q_{k-\bb{\mathbf{c}}},
\end{cases} 
&\qquad
&\begin{array}{l}
\text{if $k\in \bb{\mathbf{c}}+1$;} \\
\text{if $k\in [\bb{\mathbf{c}}+1,\bb{\mathbf{e}+1}]$,}
\end{array}
\\
\mathfrak{r}_{k} &=
\begin{cases}
\mathfrak{p}_{k},\\
\mathfrak{q}_{k-\bb{\mathbf{c}}},
\end{cases}
&\qquad
&\begin{array}{l}
\text{if $k\in \bb{\mathbf{c}}$;} \\
\text{if $k\in [\bb{\mathbf{c}},\bb{\mathbf{e}}]$,}
\end{array}
\\
V_{k} &=
\begin{cases}
T_{k},\\
U_{k-\bb{\mathbf{c}}},
\end{cases} 
&\qquad
&\begin{array}{l}
\text{if $k\in \bb{\mathbf{c}}$;} \\
\text{if $k\in [\bb{\mathbf{c}},\bb{\mathbf{e}}]$.}
\end{array}
\end{alignat*}
\end{restatable}

\begin{restatable}{proposition}{PPthComp}
\label{PPthComp}
Let $s$ be a sort in $S$, $\mathbf{c}, \mathbf{d}$ words in $S^{\star}$, $\mathfrak{P}$ a  $\mathbf{c}$-path in $\boldsymbol{\mathcal{A}}$, and $\mathfrak{Q}$ a  $\mathbf{d}$-path in $\boldsymbol{\mathcal{A}}$. Then, when defined, $\mathfrak{Q}\circ^{0}_{s}\mathfrak{P}$, the $0$-composite of $\mathfrak{P}$ and $\mathfrak{Q}$, is a $(\mathbf{c}\curlywedge\mathbf{d})$-path in $\boldsymbol{\mathcal{A}}$ of the form
$$
\xymatrix@C=50pt{
\mathfrak{Q}\circ^{0}_{s}\mathfrak{P}
\colon
\mathrm{sc}^{(0,1)}_{s}(\mathfrak{P})
\ar@{->}[r]^-{}
&
\mathrm{tg}^{(0,1)}_{s}(\mathfrak{Q}).
}
$$
Moreover, the above partial operation of $0$-composition, when defined, is associative, and, for every sort $s\in S$ and every term $P$ in $\mathrm{T}_{\Sigma}(X)_{s}$, the $(1,0)$-identity path on $P$ is, when defined, a neutral element for the operation of $0$-composition.
\end{restatable}

\begin{remark}
Later on, after defining a suitable equivalence relation on $\mathrm{Pth}_{\boldsymbol{\mathcal{A}}}$, we will show that the corresponding quotient is an $S$-sorted categorial $\Sigma$-algebra.
\end{remark}

We next define the notion of subpath of a path. But before doing that, since it will be used immediately below, we recall, from Definition~\ref{DSubw}, that, for a word $\mathbf{c}\in S^{\star}$ and indexes $k$, $l\in\bb{\mathbf{c}}$ such that $k\leq l$, $\mathbf{c}^{k,l}$ is the subword of $\mathbf{c}$ beginning at position $k$ and ending at position $l$.

\begin{restatable}{definition}{DPthSub}
\label{DPthSub}
Let $s$ be a sort in $S$, $\mathbf{c}$ a word in $S^{\star}$, $k$ and $l\in\bb{\mathbf{c}}$ such that $k\leq l$, and $\mathfrak{P}$ a $\mathbf{c}$-path in $\mathrm{Pth}_{\boldsymbol{\mathcal{A}},s}$ of the form
$$\mathfrak{P}=\left(
(P_{i})_{i\in\bb{\mathbf{c}}+1},
(\mathfrak{p}_{i})_{i\in\bb{\mathbf{c}}},
(T_{i})_{i\in\bb{\mathbf{c}}}
\right).
$$
Then we will denote by $\mathfrak{P}^{k,l}$ the ordered triple
$$
\mathfrak{P}^{k,l}
=
\left(
(P_{i}^{k,l})_{i\in\bb{\mathbf{c}^{k,l}}+1},
(\mathfrak{p}^{k,l}_{i})_{i\in\bb{\mathbf{c}^{k,l}}},
(T^{k,l}_{i})_{i\in\bb{\mathbf{c}^{k,l}}}
\right)
$$
where, 
\begin{enumerate}
\item for every $i\in \bb{\mathbf{c}^{k,l}}+1$,
\begin{multicols}{2}
\noindent
\begin{enumerate}
\item[(i)] $P_{i}^{k,l}=P_{k+i}$;
\end{enumerate}
\end{multicols}
\item for every  $i\in \bb{\mathbf{c}^{k,l}}$,
\begin{multicols}{3}
\begin{enumerate}
\item[(i)]  $\mathfrak{p}^{k,l}_{i}=\mathfrak{p}_{k+i}$;
\item[(ii)] $T^{k,l}_{i}=T_{k+i}$.
\end{enumerate}
\end{multicols}
\end{enumerate}
We will call $\mathfrak{P}^{k,l}$ the \emph{subpath of $\mathfrak{P}$ beginning at position $k$ and ending at position $l+1$}\index{path!first-order!subpath}. In particular, subpaths of the form  $\mathfrak{P}^{0,k}$ will be called \emph{initial subpaths of $\mathfrak{P}$}, and subpaths of the form  $\mathfrak{P}^{l,\bb{\mathbf{c}}-1}$ will be called \emph{final subpaths of $\mathfrak{P}$}.
\end{restatable}

In the following proposition we prove that the subpaths of a path are actually paths.

\begin{restatable}{proposition}{PPthSub}
\label{PPthSub}
Let $s$ be a sort in $S$, $\mathbf{c}$ be a word in $S^{\star}$, $k$ and $l$ indices in $\bb{\mathbf{c}}$ such that $k\leq l$, and $\mathfrak{P}$ a  $\mathbf{c}$-path in $\mathrm{Pth}_{\boldsymbol{\mathcal{A}},s}$ of the form
$$\mathfrak{P}=\left(
(P_{i})_{i\in\bb{\mathbf{c}}+1},
(\mathfrak{p}_{i})_{i\in\bb{\mathbf{c}}},
(T_{i})_{i\in\bb{\mathbf{c}}}
\right).
$$

Then $\mathfrak{P}^{k,l}$ 
is a  $\mathbf{c}^{k,l}$-path in $\mathrm{Pth}_{\boldsymbol{\mathcal{A}},s}$ of the form
$$
\xymatrix@C=50pt{
\mathfrak{P}^{k,l}
\colon
P_{k}
\ar@{->}[r]^-{}
&
P_{l+1}.
}
$$
\end{restatable}
\begin{proof}
We easily check that $P^{k,l}_{0}=P_{k+0}=P_{k}$ and $P^{k,l}_{\bb{\mathbf{c}^{k,l}}}=P_{k+l-k+1}=P_{l+1}$. Moreover, since $\mathfrak{P}$ is a $\mathbf{c}$-path in $\boldsymbol{\mathcal{A}}$ from $P_{0}$ to $P_{\bb{\mathbf{c}}}$, $\mathbf{c}^{k,l}$ is a subword of $\mathbf{c}$ and, in addition, for every $i\in \bb{\mathbf{c}^{k,l}}$, if $\mathfrak{p}^{k,l}_{i}=\mathfrak{p}_{k+i}=(M_{k+i}, N_{k+i})$, then $T_{k+i}(M_{k+i})=P_{k+i}$ and $T_{k+i}(N_{k+i})=P_{k+i+1}$. Thus, $\mathfrak{P}^{k,l}$  is a $\mathbf{c}^{k,l}$-path in $\boldsymbol{\mathcal{A}}$ from $P_{k}$ to $P_{l+1}$.
\end{proof}

We next describe how the notion of subpath of a path relates to that of $0$-composition of paths.

\begin{restatable}{proposition}{PPthRecons}
\label{PPthRecons}
Let $\mathbf{c}$ be a word in $S^{\star}$, $k$ and $l$ indexes in $\bb{\mathbf{c}}$ such that $k\leq l$, and $\mathfrak{P}$
a $\mathbf{c}$-path in $\boldsymbol{\mathcal{A}}$ of sort $s\in S$. Then
$$\mathfrak{P} = \mathfrak{P}^{l+1,\bb{\mathbf{c}}-1}\circ^{0}_{s} \mathfrak{P}^{k,l}\circ^{0}_{s}\mathfrak{P}^{0,k-1}.$$
\end{restatable}

\section{
\texorpdfstring
{Basic results on paths}
{Basic results}
}

In this section we explain how a path carries out the transformation of a term and state the structural results concerning this process. In addition, we introduce the notion of echelon, a key concept in the development of our theory.

\begin{restatable}{definition}{DEch}
\label{DEch}
We denote by $\mathrm{ech}^{(1,\mathcal{A})}$ the $S$-sorted mapping
$$
\mathrm{ech}^{(1,\mathcal{A})}
\colon
\mathcal{A}
\mor
\mathrm{Pth}_{\boldsymbol{\mathcal{A}}}
$$
which, for every sort $s$ in $S$, is defined as follows:
$$\textstyle
  \mathrm{ech}^{(1,\mathcal{A})}_{s}\nfunction
  {\mathcal{A}_{s}}
  {\mathrm{Pth}_{\boldsymbol{\mathcal{A}},s}}
  {\mathfrak{p}=(M,N)}
  {\left((M,N),\mathfrak{p},\mathrm{id}^{\mathrm{T}_{\Sigma}(X)_{s}}\right)}
$$
We also let
$
\xymatrix@C=50pt{
\mathrm{ech}^{(1,\mathcal{A})}_{s}(\mathfrak{p}) 
\colon
M
\ar@{->}[r]^-{\text{\Small{$(\mathfrak{p}, \mathrm{id}^{\mathrm{T}_{\Sigma}(X)_{s}})$}}}
&
N
}
$
stand for $\mathrm{ech}^{(1,\mathcal{A})}_{s}(\mathfrak{p})$.
This definition is sound because (i) $\mathrm{id}^{\mathrm{T}_{\Sigma}(X)_{s}}(M)=M$ and (ii) $\mathrm{id}^{\mathrm{T}_{\Sigma}(X)_{s}}(N)=N$. 
We will call $\mathrm{ech}^{(1,\mathcal{A})}_{s}(\mathfrak{p})$ the \emph{echelon associated to 
$\mathfrak{p}$}\index{echelon}\index{echelon!first-order!$\mathrm{ech}^{(1,\mathcal{A})}$}. Moreover, 
for a sort $s$ in $S$, we will say that a path $\mathfrak{P}\in \mathrm{Pth}_{\boldsymbol{\mathcal{A}},s}$ is an \emph{echelon} if there exists a rewrite rule $\mathfrak{p}\in\mathcal{A}_{s}$ such that $\mathrm{ech}^{(1,\mathcal{A})}_{s}(\mathfrak{p}) =\mathfrak{P}$, i.e., if $\mathfrak{P}\in \mathrm{Im}(\mathrm{ech}^{(1,\mathcal{A})})_{s}$. We denote by $\mathrm{Ech}^{(1,\mathcal{A})}[\mathcal{A}]$ the $S$-sorted set 
$\mathrm{Im}(\mathrm{ech}^{(1,\mathcal{A})})$. Finally, we will say that a path $\mathfrak{P}$ is \emph{echelonless}\index{path!first-order!echelonless} if $\bb{\mathfrak{P}}\geq 1$ and none of its one-step subpaths is an echelon. 
\end{restatable}

From the above it follows, as stated in the following corollary, that the translations of an echelonless path must be non-identity translations.

\begin{restatable}{corollary}{CEch}
\label{CEch}
Let $s$ be a sort in $S$, $\mathbf{c}\in S^{\star}$ and $\mathfrak{P}$ an echelonless  path in $\mathrm{Pth}_{\boldsymbol{\mathcal{A}},s}$ of the form
$$
\mathfrak{P}
=
\left(
(
P_{i}
)_{i\in\bb{\mathbf{c}}+1}, 
(
(
\mathfrak{p}_{i}
)_{i\in\bb{\mathbf{c}}},
(
T_{i}
)_{i\in\bb{\mathbf{c}}}
\right).
$$
Then $\mathbf{c}\neq \lambda$ and, for every $i\in \bb{\mathbf{c}}$, $T_{i}\neq \mathrm{id}^{\mathrm{T}_{\Sigma}(X)_{s}}$.
\end{restatable}

In the following definition we introduce the notion of a head-constant echelonless path.

\begin{restatable}{definition}{DPthHeadCt}
\label{DPthHeadCt}
Let $s$ be a sort in $S$, $\mathbf{c}$ a word in $S^{\star}-\{\lambda\}$ and  $\mathfrak{P}$ an echelonless $\mathbf{c}$-path in $\mathrm{Pth}_{\mathcal{A},s}$ of the form 
$$
\mathfrak{P}
=
\left(
(
P_{i}
)_{i\in\bb{\mathbf{c}}+1}, 
(
(
\mathfrak{p}_{i}
)_{i\in\bb{\mathbf{c}}},
(
T_{i}
)_{i\in\bb{\mathbf{c}}}
\right).
$$
We will say that $\mathfrak{P}$ is a \emph{head-constant} echelonless  path\index{path!first-order!head-constant}\index{head-constant} if $(T_{i})_{i\in\bb{\mathbf{c}}}$, the family of translations occurring in it, have the same type, i.e., they are associated to the same operation symbol. 
\end{restatable}

In the following lemma we prove that every echelonless path is head-constant.

\begin{restatable}{lemma}{LPthHeadCt}
\label{LPthHeadCt}
Let $s$ be a sort in $S$ and  $\mathfrak{P}$ an echelonless path in $\mathrm{Pth}_{\mathcal{A},s}$. Then 
$\mathfrak{P}$ is head-constant.
\end{restatable}
\begin{proof}
Let us assume that the echelonless path $\mathfrak{P}$ has the form 
$$\mathfrak{P}=
\left(
(P_{i})_{i\in\bb{\mathbf{c}}+1},
(\mathfrak{p}_{i})_{i\in\bb{\mathbf{c}}},
(T_{i})_{i\in\bb{\mathbf{c}}}
\right), 
$$
for a non-empty word $\mathbf{c}\in S^{\star}-\{\lambda\}$.

Let $i$ be an element of $\bb{\mathbf{c}}$ and let us consider the $i$-th translation occurring in 
$\mathfrak{P}$, i.e., $T_{i}$. Since $\mathfrak{P}$ is an echelonless path we have, by 
Definition~\ref{DEch}, that $T_{i}\neq \mathrm{id}^{\mathrm{T}_{\Sigma}(X)_{s}}$. Thus, for a unique word 
$\mathbf{s}_{i}\in S^{\star}-\{\lambda\}$ and a unique operation symbol $\sigma_{i}\in \Sigma_{\mathbf{s}_{i},s}$, $T_{i}$ is a translation of type $\sigma_{i}$. That is, there exists a unique index $k_{i}\in \bb{\mathbf{s}_{i}}$, a unique family of terms $(P_{i,j})_{j\in k_{i}}\in \prod_{j\in k_{i}} \mathrm{T}_{\Sigma}(X)_{s_{j}}$, a unique family of terms $(P_{i,l})_{l\in \bb{\mathbf{s}_{i}}-(k_{i}+1)}\in \prod_{l\in \bb{\mathbf{s}_{i}}-(k_{i}+1)} \mathrm{T}_{\Sigma}(X)_{s_{l}}$ and a unique translation $T'_{i}\in \mathrm{Tl}_{c_{i}}(\mathbf{T}_{\Sigma}(X))_{(s_{i})_{k_{i}}}$, the maximal prefix of $T_{i}$, see Definition~\ref{DTrans}, such that
\[
T_{i}=\sigma_{i}^{\mathbf{T}_{\Sigma}(X)}\left(
P_{i,0},\cdots, P_{i,k_{i}-1},
T'_{i}, P_{i,k_{i}+1},\cdots, P_{i,\bb{\mathbf{s}_{i}}-1}\right).
\]

Let $i$ be an element of $\bb{\mathbf{c}}-1$ and let us assume that $\mathfrak{p}_{i}=(M_{i},N_{i})$ and that $\mathfrak{p}_{i+1}=(M_{i+1},N_{i+1})$. Then, according to Definition~\ref{DPth}, the following equalities must hold
\begin{flushleft}
$\sigma_{i}^{\mathbf{T}_{\Sigma}(X)}\left(
P_{i,0},\cdots, P_{i,k_{i}-1},
T'_{i}(N_{i}), P_{i,k_{i}+1},\cdots, P_{i,\bb{\mathbf{s}_{i}}-1}\right)$
\allowdisplaybreaks
\begin{align*}
&=T_{i}(N_{i})
\\&=
P_{i+1}
\\&=T_{i+1}(M_{i+1})
\\&=
\sigma_{i+1}^{\mathbf{T}_{\Sigma}(X)}\left(
P_{i+1,0},\cdots, P_{i+1,k_{i+1}-1},
T'_{i+1}(M_{i+1}), P_{i+1,k_{i+1}+1},\cdots, P_{i+1,\bb{\mathbf{s}_{i+1}}-1}\right), 
\end{align*}
\end{flushleft}
where  $T'_{i+1}$ is the maximal prefix of $T_{i+1}$.
From this we conclude that, for every $i\in \bb{\mathbf{c}}-1$, $\mathbf{s}_{i}=\mathbf{s}_{i+1}$ and $\sigma_{i}=\sigma_{i+1}$. Thus, $(T_{i})_{i\in\bb{\mathbf{c}}}$ is a family of translations of the same type, i.e., $\mathfrak{P}$ is head-constant. 

This finishes the proof of Lemma~\ref{LPthHeadCt}.
\end{proof}

The following lemma is fundamental for the development of the present work. In it we show that every echelonless path $\mathfrak{P}$ of sort $s$ associated to an operation symbol $\sigma\in \Sigma_{\mathbf{s},s}$ has canonically associated with it a family $(\mathfrak{P}_{j})_{j\in\bb{\mathbf{s}}}$ of paths such that, for every $j\in \bb{\mathbf{s}}$, $\mathfrak{P}_{j}$ is a path of sort $s_{j}$ and they  are in charge of transforming in parallel each of the components of $\mathfrak{P}$.


\begin{restatable}{lemma}{LPthExtract}
\label{LPthExtract}\index{path!first-order!extraction procedure}
Let $s$ be a sort in $S$, $\mathbf{c}$ a word in $S^{\star}-\{\lambda\}$ and $\mathfrak{P}$
an echelonless  $\mathbf{c}$-path in $\mathrm{Pth}_{\mathcal{A},s}$ of the form 
$$\mathfrak{P}=\left(
(P_{i})_{i\in\bb{\mathbf{c}}+1},
(\mathfrak{p}_{i})_{i\in\bb{\mathbf{c}}},
(T_{i})_{i\in\bb{\mathbf{c}}}
\right).
$$
Let $\mathbf{s}$ be the unique word in $S^{\star}$ and $\sigma$ the unique operation symbol in 
$\Sigma_{\mathbf{s},s}$ for which, in virtue of Lemma~\ref{LPthHeadCt}, each of the translations of the family $(T_{i})_{i\in\bb{\mathbf{c}}}$ is of type $\sigma$. Then there exists a unique pair
$
\textstyle
((\mathbf{c}_{j})_{j\in\bb{\mathbf{s}}},(\mathfrak{P}_{j})_{j\in\bb{\mathbf{s}}})\in (S^{\star})^{\bb{\mathbf{s}}}\times \mathrm{Pth}_{\boldsymbol{\mathcal{A}},\mathbf{s}}
$
such that, for every $j\in\bb{\mathbf{s}}$, $\mathfrak{P}_{j}$ is a $\mathbf{c}_{j}$-path in $\mathrm{Pth}_{\mathcal{A},s_{j}}$ of the form
$$
\mathfrak{P}_{j}=
\left(
(P_{j,k})_{k\in\bb{\mathbf{c}_{j}}+1},
(\mathfrak{p}_{j,k})_{k\in\bb{\mathbf{c}_{j}}},
(T_{j,k})_{k\in\bb{\mathbf{c}_{j}}}
\right)$$
and there exists a unique mapping $(i^{j,k})_{(j,k)\in \coprod_{j\in \bb{\mathbf{s}}}\bb{\mathbf{c}_{j}}}\colon \coprod_{j\in \bb{\mathbf{s}}}\bb{\mathbf{c}_{j}}\mor \bb{\mathbf{c}}$ such that, for every $(j,k)$ in $\coprod_{j\in \bb{\mathbf{s}}}\bb{\mathbf{c}_{j}}$, the following equalities 
\begin{enumerate}
\item[(i)]  $c_{j,k}=c_{i^{j,k}}$ \text{ and}
\item[(ii)] $\mathfrak{p}_{j,k}=\mathfrak{p}_{i^{j,k}}$ 
\end{enumerate}
are fulfilled.
\end{restatable}

\begin{proof}
Let us assume that, for every $i\in\bb{\mathbf{c}}$, the rewrite rule $\mathfrak{p}_{i}$ is given by 
$(M_{i}, N_{i})$. Since $\mathfrak{P}$ is a path, for every $i\in\bb{\mathbf{c}}$, we have that 
\begin{enumerate}
\item[(i)] 
$T_{i}(M_{i})=P_{i}$ \text{ and}
\item[(ii)] 
$T_{i}(N_{i})=P_{i+1}$.
\end{enumerate}
Let us note that, since $\mathfrak{P}$ is echelonless, for every $i\in\bb{\mathbf{c}}$, we have that $T_{i}\neq \mathrm{id}^{\mathbf{T}_{\Sigma}(X)_{s}}$. Since each of the translations of the family 
$(T_{i})_{i\in\bb{\mathbf{c}}}$ has type $\sigma\in \Sigma_{\mathbf{s},s}$, we have that, for every $i\in\bb{\mathbf{c}}$, there exists a unique index $j_{i}\in \bb{\mathbf{s}}$, a unique family of terms $(P_{i,k})_{k\in j_{i}}\in \prod_{k\in j_{i}}\mathrm{T}_{\Sigma}(X)_{s_{k}}$, a unique family of terms $(P_{i,l})_{l\in \bb{\mathbf{s}}-(j_{i}+1)}\in \prod_{l\in \bb{\mathbf{s}}-(j_{i}+1)}\mathrm{T}_{\Sigma}(X)_{s_{l}}$ and a unique  translation $T'_{i}\in \mathrm{Tl}_{c_{i}}(\mathbf{T}_{\Sigma}(X))_{s_{j_{i}}}$  such that 
\[
T_{i}=\sigma^{\mathbf{T}_{\Sigma}(X)}\left(
P_{i,0},
\cdots,
P_{i,j_{i}-1},
T'_{i},
P_{i,j_{i}+1},
\cdots,
P_{i,\bb{\mathbf{s}}-1}
\right).
\]
In this case we will say that $\mathfrak{P}^{i,i}$ is a one-step  $c_{i}$-path of type $j_{i}$.

For every $j\in\bb{\mathbf{s}}$, let $\mathbf{c}_{j}$ be the word in $S^{\star}$ obtained by concatenating,  in order of appearance, the letters $c_{i}\in S$ for which $i\in\bb{\mathbf{c}}$ is such that $\mathfrak{P}^{i,i}$ is a one-step  $c_{i}$-path of type $j$.  

If $\mathbf{c}_{j}=\lambda$, i.e., if there are no steps of type $j$, then we define
\[\mathfrak{P}_{j}=\mathrm{ip}^{(1,0)\sharp}_{s_{j}}\left(
P_{0,j}
\right).
\]

Let us note that if $\mathfrak{P}$ does not have any step of type $j$, then the choice of the specific step does not matter, since, for every $i\in\bb{\mathbf{c}}$, we have that
$
P_{0,j}
=
P_{i,j}.
$

We will therefore focus our attention on the case in which $\mathbf{c}_{j}\neq \lambda$. 
By construction, for every $j\in\bb{\mathbf{s}}$ and every $k\in\bb{\mathbf{c}_{j}}$, there exists a unique $i^{j,k}\in\bb{\mathbf{c}}$ such that $\mathfrak{P}^{i^{j,k},i^{j,k}}$ is the $k$-th one-step path of type $j$, which starts at position $i^{j,k}$ and ends at position $i^{j,k}+1$. Hence, $c_{j,k}=c_{i^{j,k}}$, i.e., the $k$-th letter in the word $\mathbf{c}_{j}$ is equal to the letter $c_{i^{j,k}}$ in the word 
$\mathbf{c}$.

Now, for every $j\in\bb{\mathbf{s}}$, let $\mathfrak{P}_{j}$ be the ordered triple  
$$
\mathfrak{P}_{j} = \left((P_{j,k})_{k\in\bb{\mathbf{c}_{j}}+1},  (\mathfrak{p}_{j, k})_{k\in\bb{\mathbf{c}_{j}}},
(T_{ j, k})_{k\in\bb{\mathbf{c}_{j}}}
\right)
$$ 
defined as follows:
\allowdisplaybreaks
\begin{alignat*}{2}
P_{j,k} &=
\begin{cases}
T'_{i^{j,k}}\left(
M_{i^{j,k}}
\right)
, \\
T'_{i^{j,\bb{\mathbf{c}_{j}}-1}}\left(
N_{{i^{j,\bb{\mathbf{c}_{j}}-1}}}
\right),    
\end{cases} 
&\qquad
&\begin{array}{l}
\text{if $k\in \bb{\mathbf{c}_{j}}$;} \\
\text{if $k=\bb{\mathbf{c}_{j}}$;}
\end{array}
\\
\mathfrak{p}_{j,k} &=
\begin{cases}
\mathfrak{p}_{i^{j,k}},   
\end{cases} 
&\qquad
&\begin{array}{l}
\text{if $k\in \bb{\mathbf{c}_{j}}$;} 
\end{array}
\\
T'_{j,k} &=
\begin{cases}
T'_{i^{j,k}},   
\end{cases} 
&\qquad
&\begin{array}{l}
\text{if $k\in \bb{\mathbf{c}_{j}}$.} 
\end{array}
\end{alignat*}
That is, (1) for every $k\in\bb{\mathbf{c}_{j}}$, $P_{j,k}$, i.e., the $k$-th term of the sequence 
$(P_{j,k})_{k\in\bb{\mathbf{c}_{j}}+1}$, is equal to the $j$-th subterm of the term 
$T_{i^{j,k}}(M_{i^{j,k}})$, i.e., $T'_{i^{j,k}}(M_{i^{j,k}})$; (2) for every $k\in\bb{\mathbf{c}_{j}}$, $\mathfrak{p}_{j,k}$ is equal to $\mathfrak{p}_{i^{j,k}}$, i.e., the $k$-th rewrite rule of the sequence $(\mathfrak{p}_{j,k})_{k\in\bb{\mathbf{c}_{j}}}$ is equal to the rewrite rule associated to the $k$-th index of type $j$; and, finally, (3) for every $k\in\bb{\mathbf{c}_{j}}$, $T_{ j, k}$ is equal to $T'_{i^{j,k}}$, i.e., the $k$-th translation of the sequence $(T_{ j, k})_{k\in\bb{\mathbf{c}_{j}}}$ is equal to $T'_{i^{j,k}}$, the maximal prefix of the $k$-th translation of type $j$.

We claim that, for every $j\in\bb{\mathbf{s}}$, $\mathfrak{P}_{j}$ is a $\mathbf{c}_{j}$-path of sort $s_{j}$ in $\boldsymbol{\mathcal{A}}$ from 
$
T'_{i^{j,0}}(
M_{i^{j,0}}
)$ to 
$
T'_{i^{j,\bb{\mathbf{c}_{j}}-1}}(
N_{{i^{j,\bb{\mathbf{c}_{j}}-1}}}
)$. These are, respectively, the first and last term of the sequence 
$(P_{j,k})_{k\in\bb{\mathbf{c}_{j}}+1}$.

Let us note that, for every $k\in\bb{\mathbf{c}_{j}}$, we have that $\mathfrak{p}_{j, k}=\mathfrak{p}_{i^{j,k}}$, where $\mathfrak{p}_{i^{j,k}}=(M_{i^{j,k}}, N_{i^{j,k}})$, and $T_{j,k}=T'_{i^{j,k}}$. Moreover, for every $k\in\bb{\mathbf{c}_{j}}$, the following equality
\begin{enumerate}
\item[(i)] $T'_{i^{j,k}}\left(
M_{i^{j,k}}
\right)
=P_{j,k}
$
\end{enumerate}
follows, immediately, from the definition of $P_{j,k}$.

It is left to us to show that, for every $k\in\bb{\mathbf{c}_{j}}$, the following equality 
\begin{enumerate}
\item[(ii)] $T'_{i^{j,k}}\left(
N_{i^{j,k}}
\right)
=P_{j,k+1}
$
\end{enumerate}
also holds.

Let us note that for $k=\bb{\mathbf{c}_{j}}-1$, the condition (ii) above follows, immediately, by the definition of $P_{j,\bb{\mathbf{c}_{j}}}$. Therefore, it suffices to prove equality (ii) above for every $k\in\bb{\mathbf{c}_{j}}-1$. That is, we need to prove that, for every $k\in\bb{\mathbf{c}_{j}}-1$, the following equality
\[
T'_{i^{j,k}}\left(
N_{i^{j,k}}
\right)
=
T'_{i^{j,k+1}}\left(
M_{i^{j,k+1}}
\right).
\]
holds. 

Let $k\in\bb{\mathbf{c}_{j}}-1$, since $\mathfrak{P}^{i^{j,k},i^{j,k}}$, by Proposition~\ref{PPthSub}, is a $c_{i^{j,k}}$-path in $\boldsymbol{\mathcal{A}}$ from $P_{i^{j,k}}$ to $P_{i^{j,k}+1}$, we have that
\begin{itemize}
\item[(i)] $T_{i^{j,k}}\left(
M_{i^{j,k}}
\right)=P_{i^{j,k}}$ \text{ and}
\item[(ii)] $T_{i^{j,k}}\left(
N_{i^{j,k}}
\right)=P_{i^{j,k}+1}$.
\end{itemize}

Let us recall that since $T_{i^{j,k}}$ is a translation of type $j$, it has the form
\[
T_{i^{j,k}}=\sigma^{\mathbf{T}_{\Sigma}(X)}\left(
P_{i^{j,k},0},
\cdots,
P_{i^{j,k},j-1},
T'_{i^{j,k}},
P_{i^{j,k},j+1},
\cdots,
P_{i^{j,k},\bb{\mathbf{s}}-1}
\right).
\]
Moreover, for the indices after $i^{j,k}$, the corresponding derived translation would not occur at position $j$ until we reached index $i^{j,k+1}$, which would, again, be a translation of type $j$, i.e., a  translation of the form
\[
T_{i^{j,k+1}}=\sigma^{\mathbf{T}_{\Sigma}(X)}\left(
P_{i^{j,k+1},0},
\cdots,
P_{i^{j,k+1},j-1},
T'_{i^{j,k+1}},
P_{i^{j,k+1},j+1},
\cdots,
P_{i^{j,k+1},\bb{\mathbf{s}}-1}
\right).
\]

Taking this into account we would have that, for every index $i\in \bb{\mathbf{c}}$ with $i^{j,k}<i<i^{j,k+1}$, the translation $T_{i}$ would be a translation of type $j_{i}$, with $j_{i}\neq j$, i.e., a translation of the form
\[
T_{i}=\sigma^{\mathbf{T}_{\Sigma}(X)}\left(
P_{i,0},
\cdots,
P_{i,j_{i}-1},
T'_{i},
P_{i,j_{i}+1},
\cdots,
P_{i,\bb{\mathbf{s}}-1}
\right).
\]

Since $\mathfrak{P}$ is a path, the following equalities, explicitly located at position $j$,
\[
T'_{i^{j,k}}\left(
N_{i^{j,k}}
\right)
=
P_{i^{j,k}+1,j}
=
\dots
=
P_{i^{j,k+1}-1,j}
=
T'_{i^{j,k+1}}\left(
M_{i^{j,k+1}}
\right)
=
P_{j,k+1}.
\]
must be fulfilled.

This completes the proof.
\end{proof}

We emphasize that the result stated in Lemma~\ref{LPthExtract} is fundamental for the development of our theory. We will refer to it as the \emph{path extraction algorithm}. To improve the understanding of this lemma we provide the following example.

\begin{example}\label{EPthExtract}
We represent in Figure~\ref{FPthExtract} an echelonless path $\mathfrak{P}$ in $\mathrm{Pth}_{\mathbf{c},\boldsymbol{\mathcal{A}},s}$, for some $s\in S$ and $\mathbf{c}\in S^{\star}$, between complex terms. To obtain an  optimal presentation, the terms of the path have been depicted vertically. In the leftmost column we have the original path transforming a term $P_{0}$ into a term $P_{\bb{\mathbf{c}}}$, where
\begin{itemize}
\item[] $P_{0}=\sigma^{\mathbf{T}_{\Sigma}(X)}\left((P_{0, j})_{j\in\bb{\mathbf{s}}}\right)$ \text{ and}
\item[] $P_ {\bb{\mathbf{c}}}=\sigma^{\mathbf{T}_{\Sigma}(X)}\left((P_{\bb{\mathbf{c}}, j})_{j\in\bb{\mathbf{s}}}\right)$,
\end{itemize}
for a unique word $\mathbf{s}\in S^{\star}$, a unique operation symbol $\sigma\in\Sigma_{\mathbf{s},s}$, and a unique pair of families of terms $(P_{0, j})_{j\in\bb{\mathbf{s}}}$, $(P_{\bb{\mathbf{s}}, j})_{j\in\bb{\mathbf{s}}}$ in $\mathrm{T}_{\Sigma}(X)_{\mathbf{s}}$. Let us recall that, by Lemma~\ref{LPthHeadCt}, the path $\mathfrak{P}$ must be a head-constant path 
\begin{center}
\begin{tikzpicture}
[ACliment/.style={-{To [angle'=45, length=5.75pt, width=4pt, round]},
font=\scriptsize}]
\node[] (1) at (-1,0) [] {$\mathfrak{P}\colon \sigma^{\mathbf{T}_{\Sigma}(X)}((P_{0, j})_{j\in\bb{\mathbf{s}}})$};
\node[] (2) at (3.5,0) [] {$\quad$};
\node[] () at (3.5,-.075) [] {\ldots};
\node[] (3) at (7.5,0) [] {$\sigma^{\mathbf{T}_{\Sigma}(X)}((P_{\bb{\mathbf{c}}, j})_{j\in\bb{\mathbf{s}}}),$};
\draw[ACliment]  (1) to node [above]
{$(\mathfrak{p}_{0}, T_{0})$} (2);
\draw[ACliment]  (2) to node [above]
{$(\mathfrak{p}_{\bb{\mathbf{c}}-1}, T_{\bb{\mathbf{c}}-1})$} (3);
\end{tikzpicture}
\end{center}
For every $i\in \bb{\mathbf{c}}+1$, the equality $P_{i}=\sigma^{\mathbf{T}_{\Sigma}(X)}((P_{i, j})_{j\in\bb{\mathbf{s}}})$ has been depicted in the $i$-th row of Figure~\ref{FPthExtract}.

We now see how, at each step $i\in\bb{\mathbf{c}}$ and in each row, the index $i$ must be an index of type $j$, for some $j\in\bb{\mathbf{s}}$, this meaning that the rewrite rule $\mathfrak{p}_{i}$ affects only to the subterm $P_{i, j}$ of $P_{i}$ or, equivalently, that the derived translation of the $i$-th translation occurs at position $j$.

In this example, the index $0$ is a step of type $j$, since the one-step path $\mathfrak{P}^{0,0}\colon P_{0}\mor P_{1}$ only alters the subterm $P_{0, j}$ of $P_{0}$, i.e., the derived translation $T'_{0}$ occurs in the $j$-th subterm of $T_{0}$. We see how all subterms of $P_{1}$ that are different from the one  at position $j$ remain unaltered after the transformation, thus being equal to those coming from $P_{0}$. This fact is depicted in the diagram with a vertical box fitting all equal subterms. Moreover, for the subterm of type $j$, we have encountered the first transformation, from  $P_{0, j}$ to $P_{1, j}$, of the path $\mathfrak{P}_{j}$.

This process can be carried out in the index $1$, which, for this particular example, is an index of type 
$\bb{\mathbf{s}}-1$. Thus, in the one-step path $\mathfrak{P}^{1,1}\colon P_{1}\mor P_{2}$, the unique subterm that differs from $P_{1}$ to $P_{2}$ is that at position $\bb{\mathbf{s}}-1$. Hence, in the path 
$\mathfrak{P}^{0,1}\colon P_{0}\mor P_{2}$, the only subterms that have been altered are those at position $j$ and those at position $\bb{\mathbf{s}}-1$. In particular, we observe that the subterms $P_{0,0}, P_{1, 0}$, and $P_{2, 0}$ are equal, thus represented as terms in the same vertical box. As before, for the subterm of type $\bb{\mathbf{s}}-1$, we have encountered the first transformation, in this case from 
$P_{0, \bb{\mathbf{s}}-1}$ to $P_{2, \bb{\mathbf{s}}-1}$, of the path $\mathfrak{P}_{\bb{\mathbf{s}}-1}$.

This process continues for every index $i\in\bb{\mathbf{c}}$. When Figure~\ref{FPthExtract} is conceived  vertically, then we recover, for every $j\in\bb{\mathbf{s}}$, the unique $\mathbf{c}_{j}$-path $\mathfrak{P}_{j}$ in $\boldsymbol{\mathcal{A}}$ from $P_{0, j}$ to $P_{\bb{\mathbf{c}}, j}$. Let us recall that, for every $j\in\bb{\mathbf{s}}$, the $j$-th column must contain $\bb{\mathbf{c}_{j}}+1$ vertical boxes.

\begin{figure}
\centering
\begin{tikzpicture}[
fletxa/.style={thick, ->}]
\tikzstyle{every state}=[very thick, draw=blue!50,fill=blue!20, inner sep=2pt,minimum size=20pt]
\node[] (0P) at (-.8,0) [] {$P_{0}$};
\node[] (0=) at (-.2,0) [] {$=$};
\node[] (0s) at (.75,0) [] {$\sigma^{\mathbf{T}_{\Sigma}(X)}($};
\node[] (00) at (2.2,0) [] {$\color{white}P_{i, \bb{\mathbf{s}}-1}$};
\node[] () at (2.2,0) [] {$P_{0, 0}$};
\node[] (0c0) at (3.4,0) [] {$, \cdots,$};
\node[] (0j) at (4.5,0) [] {$\color{white}P_{i, \bb{\mathbf{s}}-1}$};
\node[] () at (4.5,0) [] {$P_{0, j}$};
\node[] (0c1) at (5.6,0) [] {$, \cdots,$};
\node[] (0f) at (6.7,0) [] {$P_{0, \bb{\mathbf{s}}-1}$};
\node[] (1c0) at (7.6,0) [] {$)$};

\node[] (E0label) at (-1.6,-.75) [] {$\scriptstyle (\alpha, \mathfrak{p})_{0}$};
\node[] (E0) at (-.8,-.75) [] {\rotatebox{-90}{$\mor$}};
\node[] (E00) at (2.2,-.75) [] {\rotatebox{-90}{$\,=$}};
\node[] (E0j) at (4.5,-.75) [] {\rotatebox{-90}{$\mor$}};
\node[] (E0jlabel) at (5.15,-.75) [] {$\scriptstyle (\alpha, \mathfrak{p})_{j, 0}$};
\node[] (E0f) at (6.7,-.75) [] {\rotatebox{-90}{$\,=$}};

\node[] (1P) at (-.8,-1.5) [] {$P_{1}$};
\node[] (1=) at (-.2,-1.5) [] {$=$};
\node[] (1s) at (.75,-1.5) [] {$\sigma^{\mathbf{T}_{\Sigma}(X)}($};
\node[] (10) at (2.2,-1.5) [] {$P_{1, 0}$};
\node[] (1c0) at (3.4,-1.5) [] {$, \cdots,$};
\node[] (1j) at (4.5,-1.5) [] {$\color{white}P_{i, \bb{\mathbf{s}}-1}$};
\node[] () at (4.5,-1.5) [] {$P_{1, j}$};
\node[] (1c1) at (5.6,-1.5) [] {$, \cdots,$};
\node[] (1f) at (6.7,-1.5) [] {$P_{1, \bb{\mathbf{s}}-1}$};
\node[] (1c0) at (7.6,-1.5) [] {$)$};

\node[] (E1label) at (-1.6,-2.25) [] {$\scriptstyle (\alpha, \mathfrak{p})_{1}$};
\node[] (E1) at (-.8,-2.25) [] {\rotatebox{-90}{$\mor$}};
\node[] (E10) at (2.2,-2.25) [] {\rotatebox{-90}{$\,=$}};
\node[] (E1j) at (4.5,-2.25) [] {\rotatebox{-90}{$\,=$}};
\node[] (E1f) at (6.7,-2.25) [] {\rotatebox{-90}{$\mor$}};
\node[] (E1flabel) at (7.6,-2.25) [] {$\scriptstyle (\alpha, \mathfrak{p})_{\bb{\mathbf{s}}-1, 0}$};

\node[] (2P) at (-.8,-3) [] {$P_{2}$};
\node[] (2=) at (-.2,-3) [] {$=$};
\node[] (2s) at (.75,-3) [] {$\sigma^{\mathbf{T}_{\Sigma}(X)}($};
\node[] (20) at (2.2,-3) [] {$P_{2, 0}$};
\node[] (2c0) at (3.4,-3) [] {$, \cdots,$};
\node[] (2j) at (4.5,-3) [] {$P_{2, j}$};
\node[] (2c1) at (5.6,-3) [] {$, \cdots,$};
\node[] (2f) at (6.7,-3) [] {$P_{2, \bb{\mathbf{s}}-1}$};
\node[] (2c0) at (7.6,-3) [] {$)$};

\node[] (E2label) at (-1.6,-3.75) [] {$\scriptstyle (\alpha, \mathfrak{p})_{2}$};
\node[] (E2) at (-.8,-3.75) [] {\rotatebox{-90}{$\mor$}};
\node[] (E20) at (2.2,-3.75) [] {\rotatebox{-90}{$\mor$}};
\node[] (E20label) at (2.9,-3.75) [] {$\scriptstyle (\alpha, \mathfrak{p})_{0, 0}$};
\node[] (E2jAux) at (4.5,-3.75) [] {$\color{white}P_{i, \bb{\mathbf{s}}-1}$};
\node[] (E2j) at (4.5,-3.75) [] {\rotatebox{-90}{$\, =$}};
\node[] (E2fAux) at (6.5,-3.75) [] {$\color{white}P_{i, \bb{\mathbf{s}}-1}$};
\node[] (E2f) at (6.7,-3.75) [] {\rotatebox{-90}{$\,=$}};

\node[] (3P) at (-.8,-4.5) [] {$\vdots$};
\node[] (30) at (2.2,-4.5) [] {$\vdots$};
\node[] (3c0) at (3.4,-4.5) [] {$\cdots$};
\node[] (3j) at (4.5,-4.5) [] {$\vdots$};
\node[] (3c1) at (5.6,-4.5) [] {$\cdots$};
\node[] (3f) at (6.7,-4.5) [] {$\vdots$};

\node[] (E3label) at (-1.7,-5.25) [] {$\scriptstyle (\alpha, \mathfrak{p})_{i-1}$};
\node[] (E3) at (-.8,-5.25) [] {\rotatebox{-90}{$\mor$}};
\node[] (E30Aux) at (2.2,-5.25) [] {$\color{white}P_{i, \bb{\mathbf{s}}-1}$};
\node[] (E30) at (2.2,-5.25) [] {\rotatebox{-90}{$\,=$}};
\node[] (E3j) at (4.5,-5.25) [] {\rotatebox{-90}{$\mor$}};
\node[] (E3jlabel) at (5.2,-5.25) [] {$\scriptstyle (\alpha, \mathfrak{p})_{j, k_{j}}$};
\node[] (E3fAux) at (6.7,-5.25) [] {$\color{white}P_{i, \bb{\mathbf{s}}-1}$};
\node[] (E3f) at (6.7,-5.25) [] {\rotatebox{-90}{$\,=$}};

\node[] (4P) at (-.8,-6) [] {$P_{i}$};
\node[] (4=) at (-.2,-6) [] {$=$};
\node[] (4s) at (.75,-6) [] {$\sigma^{\mathbf{T}_{\Sigma}(X)}($};
\node[] (40) at (2.2,-6) [] {$\color{white}P_{i, \bb{\mathbf{s}}-1}$};
\node[] () at (2.2,-6) [] {$P_{i, 0}$};
\node[] (4c0) at (3.4,-6) [] {$, \cdots,$};
\node[] (4j) at (4.5,-6) [] {$\color{white}P_{i, \bb{\mathbf{s}}-1}$};
\node[] (4) at (4.5,-6) [] {$P_{i, j}$};
\node[] (4c1) at (5.6,-6) [] {$, \cdots,$};
\node[] (4f) at (6.7,-6) [] {$P_{i,\bb{\mathbf{s}}-1}$};
\node[] (4c0) at (7.6,-6) [] {$)$};

\node[] (E4label) at (-1.6,-6.75) [] {$\scriptstyle (\alpha, \mathfrak{p})_{i}$};
\node[] (E4) at (-.8,-6.75) [] {\rotatebox{-90}{$\mor$}};
\node[] (E40Aux) at (2.2,-6.75) [] {$\color{white}P_{i, \bb{\mathbf{s}}-1}$};
\node[] (E40) at (2.2,-6.75) [] {\rotatebox{-90}{$\,=$}};
\node[] (E4jAux) at (4.5,-6.75) [] {$\color{white}P_{i, \bb{\mathbf{s}}-1}$};
\node[] (E4j) at (4.5,-6.75) [] {\rotatebox{-90}{$\,=$}};
\node[] (E4f) at (6.7,-6.75) [] {\rotatebox{-90}{$\mor$}};
\node[] (E3jlabel) at (7.9,-6.75) [] {$\scriptstyle (\alpha, \mathfrak{p})_{\bb{\mathbf{s}}-1, k_{\bb{\mathbf{s}}-1}}$};

\node[] (5P) at (-.8,-7.5) [] {$\vdots$};
\node[] (50) at (2.2,-7.5) [] {$\vdots$};
\node[] (5c0) at (3.4,-7.5) [] {$\cdots$};
\node[] (5j) at (4.5,-7.5) [] {$\vdots$};
\node[] (5c1) at (5.6,-7.5) [] {$\cdots$};
\node[] (5f) at (6.7,-7.5) [] {$\vdots$};

\node[] (E5label) at (-1.6,-8.25) [] {$\scriptstyle (\alpha, \mathfrak{p})_{_{\bb{\mathbf{c}}-1}}$};
\node[] (E5) at (-.8,-8.25) [] {\rotatebox{-90}{$\mor$}};
\node[] (E50) at (2.2,-8.25) [] {\rotatebox{-90}{$\mor$}};
\node[] (E50label) at (2.9,-8.25) [] {$\scriptstyle (\alpha, \mathfrak{p})_{0, k_{0}}$};
\node[] (E5jAux) at (4.5,-8.25) [] {$\color{white}P_{i, \bb{\mathbf{s}}-1}$};
\node[] (E5j) at (4.5,-8.25) [] {\rotatebox{-90}{$\, =$}};
\node[] (E5fAux) at (6.7,-8.25) [] {$\color{white}P_{i, \bb{\mathbf{s}}-1}$};
\node[] (E5f) at (6.7,-8.25) [] {\rotatebox{-90}{$\,=$}};

\node[] (6P) at (-.8,-9) [] {$P_{\bb{\mathbf{c}}}$};
\node[] (6=) at (-.2,-9) [] {$=$};
\node[] (6s) at (.75,-9) [] {$\sigma^{\mathbf{T}_{\Sigma}(X)}($};
\node[] (60) at (2.2,-9) [] {$\color{white}P_{i, \bb{\mathbf{s}}-1}$};
\node[] () at (2.2,-9) [] {$\color{white}P_{i, \bb{\mathbf{s}}-1}$};
\node[] () at (2.2,-9) [] {$P_{\bb{\mathbf{c}}, 0}$};
\node[] (6c0) at (3.4,-9) [] {$, \cdots,$};
\node[] (6j) at (4.5,-9) [] {$\color{white}P_{i, \bb{\mathbf{s}}-1}$};
\node[] (6) at (4.5,-9) [] {$P_{\bb{\mathbf{c}}, j}$};
\node[] (6c1) at (5.6,-9) [] {$, \cdots,$};
\node[] (6f) at (6.7,-9) [] {$\color{white}P_{i, \bb{\mathbf{s}}-1}$};
\node[] () at (6.7,-9) [] {$P_{\bb{\mathbf{c}}, \bb{\mathbf{s}}-1}$};
\node[] (6c0) at (7.6,-9) [] {$)$};

\tikzset{encercla/.style={draw=black, line width=.5pt, inner sep=0pt, rectangle, rounded corners}};
\node [encercla,  fit=(0f) (1f) ] {} ;
\node [encercla, fit=(2f) (3f) (4f) ] {} ;
\node [encercla, fit=(00) (10) (20)] {} ;
\node [encercla, fit=(30) (40) (50)] {} ;
\node [encercla, fit=(0j)] {} ;
\node [encercla, fit= (1j) (2j) (3j)] {} ;
\node [encercla, fit=(4j) (6j)] {} ;
\node [encercla, fit=(5f) (6f)] {} ;
\node [encercla, fit=(60) ] {} ;

\tikzset{encercla2/.style={draw=white, line width=.8pt, inner sep=0pt, rectangle}};

\node [encercla2, dashed, fit=(E2jAux)] {} ;
\node [encercla2, dashed, fit=(E2fAux)] {} ;
\node [encercla2, dashed, fit=(E30Aux)] {} ;
\node [encercla2, dashed, fit=(E3fAux)] {} ;
\node [encercla2, dashed, fit=(E40Aux)] {} ;
\node [encercla2, dashed, fit=(E4jAux)] {} ;
\node [encercla2, dashed, fit=(E5fAux)] {} ;
\node [encercla2, dashed, fit=(E5jAux)] {} ;
\end{tikzpicture}
\caption{An example of the path extraction algorithm for an echelonless path.}\label{FPthExtract}
\end{figure}

On the other hand, the diagram in Figure~\ref{FPthExtractProc} provides an overall summary of Lemma~\ref{LPthExtract}. Let us recall that the just mentioned lemma does not entails that, for every $j\in\bb{\mathbf{s}}$, the extracted paths must necessarily be non-trivial. It could be the case that the overall path does not alter the subterm $P_{0, j}$, for some $j\in\bb{\mathbf{s}}$, thus leaving $\mathfrak{P}_{j}$ as the $(1,0)$-identity path on $P_{0, j}$. It could also happen that, for some $j\in\bb{\mathbf{s}}$, the extracted path $\mathfrak{P}_{j}$ contains at least one echelon, however one should note that this path will start with a term of height smaller than that of the source of the original path. Thus giving rise, as far as paths are concerned, to the possibility of making proofs by induction and definitions by recursion, as we shall see later on.

\begin{figure}
\centering
\begin{tikzpicture}[
fletxa/.style={thick, ->}]
\tikzstyle{every state}=[very thick, draw=blue!50,fill=blue!20, inner sep=2pt,minimum size=20pt]
\node[] () at (-.65,0) [] {$P_{0}$};
\node[] () at (-1.1,-.75)[] {$\mathfrak{P}$};
\node[] () at (-.65,-.75)[] {\rotatebox{-90}{$\mor$}};
\node[] () at (-.65,-1.5) [] {$P_{\bb{\mathbf{c}}}$};
\node[] () at (.475,.05) [] {$=\sigma^{\mathbf{T}_{\Sigma}(X)}($};
\node[] () at (.475,-1.45) [] {$=\sigma^{\mathbf{T}_{\Sigma}(X)}($};

\node[] () at (1.75,0) [] {$P_{0, 0}$};
\node[] () at (2.5,-.75)[] {$\mathfrak{P}_{0}\, \cdots$};
\node[] () at (1.75,-.75)[] {\rotatebox{-90}{$\mor$}};
\node[] () at (1.75,-1.5) [] {$P_{\bb{\mathbf{c}}, 0}$};

\node[] () at (2.75,0) [] {$, \cdots, $};
\node[] () at (2.75,-1.5) [] {$, \cdots,$};

\node[] () at (3.75,0) [] {$P_{0, j}$};
\node[] () at (4.5,-.75)[] {$\mathfrak{P}_{j}\,\cdots$};
\node[] () at (3.75,-.75)[] {\rotatebox{-90}{$\mor$}};
\node[] () at (3.75,-1.5) [] {$P_{\bb{\mathbf{c}}, j}$};

\node[] () at (4.75,0) [] {$, \cdots, $};
\node[] () at (4.75,-1.5) [] {$, \cdots,$};

\node[] () at (5.9,0) [] {$P_{0, \bb{\mathbf{s}}-1})$};
\node[] () at (6.6,-.75)[] {$\mathfrak{P}_{\bb{\mathbf{s}}-1}$};
\node[] () at (5.9,-.75)[] {\rotatebox{-90}{$\mor$}};
\node[] () at (5.9,-1.5) [] {$P_{\bb{\mathbf{c}}, \bb{\mathbf{s}}-1})$};
\end{tikzpicture}
\caption{Display of the path extraction algorithm.}\label{FPthExtractProc}
\end{figure}
\end{example}
\chapter{
\texorpdfstring
{Algebraic structures on $\mathrm{Pth}_{\boldsymbol{\mathcal{A}}}$}
{Algebraic structures on paths}
}\label{S1B}

In this chapter we begin by defining a structure of $\Sigma$-algebra on $\mathrm{Pth}_{\boldsymbol{\mathcal{A}}}$, the $S$-sorted set of paths in $\boldsymbol{\mathcal{A}}$, we denote by $\mathbf{Pth}^{(0,1)}_{\boldsymbol{\mathcal{A}}}$ the corresponding $\Sigma$-algebra. To do this, we define for every pair $(\mathbf{s},s)\in S^{\star}\times S$, every operation symbol $\sigma\in\Sigma_{\mathbf{s},s}$ and every family of paths $(\mathfrak{P}_{j})_{j\in\bb{\mathbf{s}}}$ in $\mathrm{Pth}_{\boldsymbol{\mathcal{A}},\mathbf{s}}$, the path 
$\sigma^{\mathbf{Pth}_{\boldsymbol{\mathcal{A}}}}((\mathfrak{P}_{j})_{j\in\bb{\mathbf{s}}})$. We then prove that 
$\sigma^{\mathbf{Pth}_{\boldsymbol{\mathcal{A}}}}((\mathfrak{P}_{j})_{j\in\bb{\mathbf{s}}})$ is always an echelonless path and that the path extraction algorithm applied to it always returns the original path family $(\mathfrak{P}_{j})_{j\in\bb{\mathbf{s}}}$. Then we prove that the one-step echelonless paths can be completely characterized by means of the path extraction algorithm and the composition defined above. In addition, we prove that the path $\sigma^{\mathbf{Pth}_{\boldsymbol{\mathcal{A}}}}((\mathfrak{P}_{j})_{j\in\bb{\mathbf{s}}})$, when performed on a family of $(1,0)$-identity paths $(\mathfrak{P}_{j})_{j\in\bb{\mathbf{s}}}$ is itself a $(1,0)$-identity path. This allows us to obtain the subalgebra of $(1,0)$-identity paths, denoted by $\mathrm{ip}^{(1,0)\sharp}[\mathbf{T}_{\Sigma}(X)]$, of $\mathbf{Pth}^{(0,1)}_{\boldsymbol{\mathcal{A}}}$. We next show that the $S$-sorted mappings $\mathrm{sc}^{(0,1)}$, $\mathrm{tg}^{(0,1)}$ and $\mathrm{ip}^{(1,0)\sharp}$, of $(0,1)$-source, $(0,1)$-target and $(1,0)$-identity paths, respectively, are $\Sigma$-homomorphisms. In fact, we show that $\mathrm{ip}^{(1,0)\sharp}$, the $(1,0)$-identity path mapping, can be obtained, by the universal property of $\mathbf{T}_{\Sigma}(X)$, as the homomorphic extension of the $S$-sorted mapping $\mathrm{ip}^{(1,X)}$, that, for every $s\in S$, maps every variable in $X_{s}$ to the $(1,0)$-identity path on it. We, finally, prove that the $S$-sorted $\Sigma$-algebras $\mathbf{Pth}^{(0,1)}_{\boldsymbol{\mathcal{A}}}/{\mathrm{Ker}(\mathrm{sc}^{(0,1)})}$, $\mathbf{Pth}^{(0,1)}_{\boldsymbol{\mathcal{A}}}/{\mathrm{Ker}(\mathrm{tg}^{(0,1)})}$ and $\mathbf{T}_{\Sigma}(X)$ are isomorphic.


\begin{restatable}{proposition}{PPthAlg}
\label{PPthAlg}\index{path!first-order!$\mathbf{Pth}^{(0,1)}_{\boldsymbol{\mathcal{A}}}$}
The $S$-sorted set $\mathrm{Pth}_{\boldsymbol{\mathcal{A}}}$ is equipped, in a natural way, with a structure of $\Sigma$-algebra.
\end{restatable}

\begin{proof}
Let us denote by $\mathbf{Pth}^{(0,1)}_{\boldsymbol{\mathcal{A}}}$ the $\Sigma$-algebra defined as follows

\textsf{(1)} The underlying $S$-sorted set of $\mathbf{Pth}^{(0,1)}_{\boldsymbol{\mathcal{A}}}$ is $\mathrm{Pth}_{\boldsymbol{\mathcal{A}}} = (\mathrm{Pth}_{\boldsymbol{\mathcal{A}},s})_{s\in S}$.

\textsf{(2)} For every $(\mathbf{s},s)\in S^{\star}\times S$ and every operation symbol $\sigma\in\Sigma_{\mathbf{s},s}$, the operation $\sigma^{\mathbf{Pth}^{(0,1)}_{\boldsymbol{\mathcal{A}}}}$, abbreviated to $\sigma^{\mathbf{Pth}_{\boldsymbol{\mathcal{A}}}}$  , from $\mathrm{Pth}_{\boldsymbol{\mathcal{A}}, \mathbf{s}}$ to $\mathrm{Pth}_{\boldsymbol{\mathcal{A}},s}$ associated to $\sigma$ assigns to a family of paths
$(\mathfrak{P}_{j})_{j\in\bb{\mathbf{s}}}\in \mathrm{Pth}_{\boldsymbol{\mathcal{A}}, \mathbf{s}}$ where, for every $j\in\bb{\mathbf{s}}$, $\mathfrak{P}_{j}$ is a $\mathbf{c}_{j}$-path in $\boldsymbol{\mathcal{A}}$ from $P_{j,0}$ to $P_{j, \bb{\mathbf{c}_{j}}}$ of the form
$$\mathfrak{P}_{j}=
\left((P_{j,k})_{k\in\bb{\mathbf{c}_{j}}+1},
(\mathfrak{p}_{j,k})_{k\in\bb{\mathbf{c}_{j}}},
(T_{j,k})_{k\in\bb{\mathbf{c}_{j}}}
\right),
$$
for a  unique $\mathbf{c}_{j}$ in $ S^{\star}$, and a unique pair of terms $P_{j,0}$,$P_{j, \bb{\mathbf{c}_{j}}}$ in $\mathrm{T}_{\Sigma}(X)_{s_{j}}$, precisely the $\mathbf{c}$-path in $\boldsymbol{\mathcal{A}}$ of sort $s$ given by
$$
\sigma^{\mathbf{Pth_{\boldsymbol{\mathcal{A}}}}}
\left(
(\mathfrak{P}_{j})_{j\in\bb{\mathbf{s}}}
\right)
=
\left(
(P_{i})_{i\in\bb{\mathbf{c}}+1},
(\mathfrak{p}_{i})_{i\in\bb{\mathbf{c}}},
(T_{i})_{i\in\bb{\mathbf{c}}},
\right),
$$
where $\mathbf{c}=\bigcurlywedge_{j\in\bb{\mathbf{s}}}\mathbf{c}_{j}$ is the concatenation of the family of words $(\mathbf{c}_{j})_{j\in\bb{\mathbf{s}}}$.

Before verifying that the definition $\sigma^{\mathbf{Pth_{\boldsymbol{\mathcal{A}}}}}((\mathfrak{P}_{j})_{j\in\bb{\mathbf{s}}})$ is sound, let us point out the following facts. By construction of $\mathbf{c}$, it happens that $\bb{\mathbf{c}}=\sum_{j\in\bb{\mathbf{s}}}\bb{\mathbf{c}_{j}}$. Hence, the $i$-th letter of $\mathbf{c}$ will be the $k$-th letter of the subword $\mathbf{c}_{j}$, for a unique $j\in\bb{\mathbf{s}}$ and a unique $k\in\bb{\mathbf{c}_{j}}$. We will write $i=(j,k)$ to denote this dependency. Moreover, for the aforementioned indices, we have that $c_{i}=c_{j, k}$.

Returning to the definition of $\sigma^{\mathbf{Pth_{\boldsymbol{\mathcal{A}}}}}((\mathfrak{P}_{j})_{j\in\bb{\mathbf{s}}})$, for $i\in\bb{\mathbf{c}}$ with $i=(j,k)$, we define the $0$-constituent at step $i$ of $\sigma^{\mathbf{Pth_{\boldsymbol{\mathcal{A}}}}}((\mathfrak{P}_{j})_{j\in\bb{\mathbf{s}}})$ to be the term
\[
P_{i}  =
\sigma^{\mathbf{T}_{\Sigma}(X)}\left(P_{0,\bb{\mathbf{c}_{0}}},\cdots,P_{j-1, \bb{\mathbf{c}_{j-1}}}, P_{j, k},P_{ j+1,0},\cdots, P_{\bb{\mathbf{s}}-1, 0}\right).
\]

That is, if $i\in\bb{\mathbf{c}}$ and $i=(j,k)$, then we have that, to the left of position $j$, every subterm is equal to the last term of the corresponding path, and, to the right of position $j$, every subterm is equal to the initial term of the corresponding path. The $j$-th subterm of $P_{i}$ is the $k$-th term appearing in the path $\mathfrak{P}_{j}$. In particular, since $0=(0,0)$, we have that 
\[P_{0}=\sigma^{\mathbf{T}_{\Sigma}(X)}\left(\left(
P_{j,0}
\right)_{j\in\bb{\mathbf{s}}}
\right).\]
Finally, for the case $i=\bb{\mathbf{c}}$, we define
\[
P_{\bb{\mathbf{c}}}  =
\sigma^{\mathbf{T}_{\Sigma}(X)}\left(\left(
P_{j,\bb{\mathbf{c}_{j}}}
\right)_{j\in\bb{\mathbf{s}}}
\right).
\]

For $i\in\bb{\mathbf{c}}$ with $i=(j,k)$, we define the rewrite rule $\mathfrak{p}_{i}$ to be equal to $\mathfrak{p}_{j, k}$. That is, the $i$-th rewrite rule of $\sigma^{\mathbf{Pth}_{\boldsymbol{\mathcal{A}}}}((\mathfrak{P}_{j})_{j\in\bb{\mathbf{s}}})$ is equal to the $k$-th rewrite rule of the path $\mathfrak{P}_{j}$. 

Finally, for $i\in\bb{\mathbf{c}}$ with $i=(j,k)$, we define the translation at step $i$  
of $\sigma^{\mathbf{Pth}_{\boldsymbol{\mathcal{A}}}}((\mathfrak{P}_{j})_{j\in\bb{\mathbf{s}}})$ to be equal to
\[
T_{i}=
\sigma^{\mathbf{T}_{\Sigma}(X)}
\left(
P_{0,\bb{\mathbf{c}_{0}}}
,
\cdots,
P_{j-1,\bb{\mathbf{c}_{j-1}}}
,
T_{j,k},
P_{j+1,0}
\cdots,
P_{\bb{\mathbf{s}}-1,0}
\right).
\]
That is, if $i\in\bb{\mathbf{c}}$ and $i=(j,k)$, then we have that, to the left of position $j$, every subterm is equal to the  term of the last $0$-constituent of the corresponding path, and, to the right of position $j$, every subterm is equal to the term of the initial $0$-constituent  of the corresponding  path.  The $j$-th subterm of $T_{i}$ is the $k$-th translation appearing in the path $\mathfrak{P}_{j}$.  Note that, by Definition~\ref{DTrans}, for every $i\in\bb{\mathbf{c}}$, $T_{i}$ is a translation.

\begin{claim}\label{CPthSigma} Let $(\mathbf{s},s)$ be a pair in $S^{\star}\times S$ and, for every $j\in\bb{\mathbf{s}}$, let $\mathfrak{P}_{j}$  be a $\mathbf{c}_{j}$-path in $\mathrm{Pth}_{\boldsymbol{\mathcal{A}},s_{j}}$. Then $\sigma^{\mathbf{Pth}_{\boldsymbol{\mathcal{A}}}}((\mathfrak{P}_{j})_{j\in\bb{\mathbf{s}}})$ is a $\bigcurlywedge_{j\in\bb{\mathbf{s}}}\mathbf{c}_{j}$-path in $\boldsymbol{\mathcal{A}}$ of the form
$$
\sigma^{\mathbf{Pth}_{\boldsymbol{\mathcal{A}}}}
\left(\left(\mathfrak{P}_{j}
\right)_{j\in\bb{\mathbf{s}}}
\right)
\colon
\sigma^{\mathbf{T}_{\Sigma}(X)}\left(\left(
\mathrm{sc}^{(0,1)}_{s_{j}}\left(
\mathfrak{P}_{j}\right)
\right)_{j\in\bb{\mathbf{s}}}\right)
\mor
\sigma^{\mathbf{T}_{\Sigma}(X)}\left(\left(
\mathrm{tg}^{(0,1)}_{s_{j}}\left(\mathfrak{P}_{j}
\right)
\right)_{j\in\bb{\mathbf{s}}}\right).
$$
\end{claim}

The first two conditions can be easily checked, since $P_{0}=\sigma^{\mathbf{T}_{\Sigma}(X)}((P_{j,0})_{j\in\bb{\mathbf{s}}})$ and $P_{\bb{\mathbf{c}}}=\sigma^{\mathbf{T}_{\Sigma}(X)}((P_{j,\bb{\mathbf{c}_{j}}})_{j\in\bb{\mathbf{s}}})$. 

Note that for every $i\in\bb{\mathbf{c}}$, there exists $j\in\bb{\mathbf{s}}$ and $k\in\bb{\mathbf{c}_{j}}$ such that $i=(j,k)$ or, equivalently, $c_{i}=c_{j,k}$. Hence, $\mathfrak{p}_{i}=\mathfrak{p}_{j,k}$, where $\mathfrak{p}_{j,k}$ is the $k$-th  rewrite rule of the $j$-th  path $\mathfrak{P}_{j}$. Let us assume that $\mathfrak{p}_{j,k}=(M_{j,k}, N_{j,k})$. Then, since $\mathfrak{P}_{j}$ is a $\mathbf{c}_{j}$-path in $\boldsymbol{\mathcal{A}}$ from $P_{j,0}$ to $P_{j,\bb{\mathbf{c}}_{j}}$, it follows that
\begin{multicols}{2}
\begin{itemize}
\item[(i)] $T_{j,k}(M_{j,k})=P_{j,k}$;
\item[(ii)] $T_{j,k}(N_{j,k})=P_{j,k+1}$.
\end{itemize}
\end{multicols}

Note that the following equality holds
\allowdisplaybreaks
\begin{align*}
T_{i}\left(
M_{i}
\right)&=
\sigma^{\mathbf{T}_{\Sigma}(X)}
\left(
P_{0,\bb{\mathbf{c}_{0}}},
\cdots,
P_{j-1,\bb{\mathbf{c}_{j-1}}}
,
T_{j,k}\left(
M_{j,k}
\right),
P_{j+1,0},
\cdots,
P_{\bb{\mathbf{s}}-1,0}
\right)
\tag{1}
\\&=
\sigma^{\mathbf{T}_{\Sigma}(X)}
\left(
P_{0,\bb{\mathbf{c}_{0}}},
\cdots,
P_{j-1,\bb{\mathbf{c}_{j-1}}}
,
P_{j,k},
P_{j+1,0},
\cdots,
P_{\bb{\mathbf{s}}-1,0}
\right)
\tag{2}
\\&=
P_{i}.
\tag{3}
\end{align*}

In the just stated chain of equalities, the first equality recovers the description of $T_{i}$ and $M_{i}$; the second equality follows from item (i) above; finally, the third equality recovers the description of $P_{i}$.

Now, for the index $i+1$, either (1) the $(i+1)$-th letter of $\mathbf{c}$ is also a letter of the subword $\mathbf{c}_{j}$, and therefore $i+1=(j,k+1)$ or (2) the $(i+1)$-th letter of $\mathbf{c}$ is the first letter of the subword $\mathbf{c}_{j+1}$, hence $i+1=(j+1,0)$ and, necessarily, $i=(j,\bb{\mathbf{c}_{j}}-1)$. These two cases can be handled similarly. Taking into account the statements above, we conclude that 
$
T_{i}\left(
N_{i}
\right)
=
P_{i+1}
.
$

Therefore $\sigma^{\mathbf{Pth}_{\boldsymbol{\mathcal{A}}}}((\mathfrak{P}_{j})_{j\in\bb{\mathbf{s}}})$ is a $\mathbf{c}$-path in $\boldsymbol{\mathcal{A}}$ from  
$\sigma^{\mathbf{T}_{\Sigma}(X)}((\mathrm{sc}^{(0,1)}_{s_{j}}(\mathfrak{P}_{j}))_{j\in\bb{\mathbf{s}}})$ to $\sigma^{\mathbf{T}_{\Sigma}(X)}(\mathrm{tg}^{(0,1)}_{s_{j}}(\mathfrak{P}_{j}))_{j\in\bb{\mathbf{s}}})$.

This finishes the proof.
\end{proof}

We conclude this section by stating some basic properties about the just defined many-sorted $\Sigma$-algebra of paths $\mathbf{Pth}_{\boldsymbol{\mathcal{A}}}^{(0,1)}$. An immediate consequence of Proposition~\ref{PPthAlg} is that every non-constant operation symbol of the signature $\Sigma$, regardles of its arity and coarity, when applied to a family of paths (not entirely composed of $(1,0)$-identity paths) always return an echelonless path.

\begin{restatable}{corollary}{CPthWB}
\label{CPthWB} Let $(\mathbf{s},s)$ be an element of $(S^{\star}-\{\lambda\})\times S$, $\sigma$ an operation symbol in $\Sigma_{\mathbf{s},s}$, and $(\mathfrak{P}_{j})_{j\in\bb{\mathbf{s}}}$ a family of paths in $\mathrm{Pth}_{\boldsymbol{\mathcal{A}},\mathbf{s}}$ such that, for some $j\in\bb{\mathbf{s}}$, $\mathfrak{P}_{j}$ is a path of length at least one. Then  the path $\sigma^{\mathbf{Pth}_{\boldsymbol{\mathcal{A}}}}((\mathfrak{P}_{j})_{j\in\bb{\mathbf{s}}})$ is an echelonless path.
\end{restatable}

We next state that the application of the path extraction algorithm, described in Lemma~\ref{LPthExtract}, to a path obtained from the action of a non-constant operation symbol on a family of paths not entirely composed of $(1,0)$-identity paths recovers its original components. In this construction, however, the appearance of the paths is no longer chaotic, since the operations introduced in Proposition~\ref{PPthAlg} follow a leftmost innermost strategy. The description of the overall process is represented in Figure~\ref{FPthSigma}.

\begin{restatable}{proposition}{PRecov}
\label{PRecov}
Let $(\mathbf{s},s)$ be an element of $ (S^{\star}-\{\lambda\})\times S$, $\sigma$ an operation symbol in $\Sigma_{\mathbf{s},s}$, and $(\mathfrak{P}_{j})_{j\in\bb{\mathbf{s}}}$ a family of paths in $\mathrm{Pth}_{\boldsymbol{\mathcal{A}},\mathbf{s}}$. If, for some $j\in\bb{\mathbf{s}}$, $\mathfrak{P}_{j}$ is a path of length at least one, then the path extraction algorithm from Lemma~\ref{LPthExtract} applied to the path $\sigma^{\mathbf{Pth}_{\boldsymbol{\mathcal{A}}}}((\mathfrak{P}_{j})_{j\in\bb{\mathbf{s}}})$ recovers the original family $(\mathfrak{P}_{j})_{j\in\bb{\mathbf{s}}}$.
\end{restatable}

\begin{figure}
\centering
\begin{scaletikzpicturetowidth}{\textwidth}
\begin{tikzpicture}[
scale=\tikzscale,
fletxa/.style={thick, ->}]
\tikzstyle{every state}=[very thick, draw=blue!50,fill=blue!20, inner sep=2pt,minimum size=20pt]
\def\posPx{-1.55}
\def\posEqx{-.6}
\def\posSigmax{.25}
\def\posZerox{2.15}
\def\posDotOnex{3.8}
\def\posJx{5.45}
\def\posDotTwox{7.15}
\def\posFx{8.8}
\def\posEndx{10.1}

\node[] (00) at (\posZerox,0) [] 
{$\color{white}P_{\bb{\mathbf{c}_{\bb{\mathbf{s}}-1}}, \bb{\mathbf{s}}-1}$};
\node[] (10) at (\posZerox,-1.5) [] 
{$\color{white}P_{\bb{\mathbf{c}_{\bb{\mathbf{s}}-1}}, \bb{\mathbf{s}}-1}$};
\node[] (20) at (\posZerox,-3) [] 
{$\color{white}P_{\bb{\mathbf{c}_{\bb{\mathbf{s}}-1}}, \bb{\mathbf{s}}-1}$};
\node[] (30) at (\posZerox,-4.4) [] {};
\node[] (40) at (\posZerox,-6) [] 
{$\color{white}P_{\bb{\mathbf{c}_{\bb{\mathbf{s}}-1}}, \bb{\mathbf{s}}-1}$};
\node[] (50) at (\posZerox,-7.4) [] {};
\node[] (60) at (\posZerox,-9) []
{$\color{white}P_{\bb{\mathbf{c}_{\bb{\mathbf{s}}-1}}, \bb{\mathbf{s}}-1}$};
\node[] (70) at (\posZerox,-10.25) [] {};
\node[] (80) at (\posZerox,-11.75) []
{$\color{white}P_{\bb{\mathbf{c}_{\bb{\mathbf{s}}-1}}, \bb{\mathbf{s}}-1}$};
\node[] (90) at (\posZerox,-13.24) [] {};
\node[] (100) at (\posZerox,-14.75) []
{$\color{white}P_{\bb{\mathbf{c}_{\bb{\mathbf{s}}-1}}, \bb{\mathbf{s}}-1}$};

\node[] (0j) at (\posJx,0) [] 
{$\color{white}P_{\bb{\mathbf{c}_{\bb{\mathbf{s}}-1}}, \bb{\mathbf{s}}-1}$};
\node[] (1j) at (\posJx,-1.4) [] {};
\node[] (2j) at (\posJx,-3) [] {};
\node[] (3j) at (\posJx,-4.4) [] {};
\node[] (4j) at (\posJx,-6) [] 
{$\color{white}P_{\bb{\mathbf{c}_{\bb{\mathbf{s}}-1}}, \bb{\mathbf{s}}-1}$};
\node[] (5j) at (\posJx,-7.5) []
{$\color{white}P_{\bb{\mathbf{c}_{\bb{\mathbf{s}}-1}}, \bb{\mathbf{s}}-1}$};
\node[] (6j) at (\posJx,-9) []
{$\color{white}P_{\bb{\mathbf{c}_{\bb{\mathbf{s}}-1}}, \bb{\mathbf{s}}-1}$};
\node[] (7j) at (\posJx,-10.25) [] {};
\node[] (8j) at (\posJx,-11.75) []
{$\color{white}P_{\bb{\mathbf{c}_{\bb{\mathbf{s}}-1}}, \bb{\mathbf{s}}-1}$};
\node[] (9j) at (\posJx,-13.24) [] {};
\node[] (10j) at (\posJx,-14.75) []
{$\color{white}P_{\bb{\mathbf{c}_{\bb{\mathbf{s}}-1}}, \bb{\mathbf{s}}-1}$};

\node[] (0f) at (\posFx,0) []
{$\color{white}P_{\bb{\mathbf{c}_{\bb{\mathbf{s}}-1}}, \bb{\mathbf{s}}-1}$};
\node[] (1f) at (\posFx,-1.4) [] {};
\node[] (2f) at (\posFx,-3) [] {};
\node[] (3f) at (\posFx,-4.4) [] {};
\node[] (4f) at (\posFx,-6) [] {};
\node[] (5f) at (\posFx,-7.5) [] {};
\node[] (6f) at (\posFx,-9) [] {};
\node[] (7f) at (\posFx,-10.25) [] {};
\node[] (8f) at (\posFx,-11.75) [] {$P_{\bb{\mathbf{s}}-1, 0}$};
\node[] (9f) at (\posFx,-13.3) []
{$\color{white}P_{\bb{\mathbf{c}_{\bb{\mathbf{s}}-1}}, \bb{\mathbf{s}}-1}$};
\node[] (10f) at (\posFx,-14.75) [] {$P_{\bb{\mathbf{s}}-1, \bb{\mathbf{c}_{\bb{\mathbf{s}}-1}}}$};

\tikzset{encercla/.style={draw=black, 
pattern=crosshatch dots, pattern color=blue!20!white, line width=.5pt, inner sep=0pt, rectangle, rounded corners}};
\tikzset{encerclad/.style={draw=black, 
pattern=dots, pattern color=blue!20!white, line width=.5pt, inner sep=0pt, rectangle, rounded corners}};
\node [encercla,  fit=(00) ] {} ;
\node [encerclad, dashed, fit=(10)] {} ;
\node [encercla, fit=(20) (30) (40) (50) (60) (70) (80) (90) (100)] {} ;
\node [encercla,  fit=(0j) (1j) (2j) (3j) (4j)] {} ;
\node [encerclad, dashed,  fit=(5j)] {} ;
\node [encercla, fit=(6j) (7j) (8j) (9j) (10j)] {} ;
\node [encercla,  fit=(0f) (1f) (2f) (3f) (4f) (5f) (6f) (7f) (8f)] {} ;
\node [encercla, fit=(10f)] {} ;
\node [encerclad, dashed,  fit=(9f)] {} ;

\node[] (0P) at (\posPx,0) [] {$P_{0}$};
\node[] (0=) at (\posEqx,0) [] {$=$};
\node[] (0s) at (\posSigmax,.05) [] {$\sigma^{\mathbf{T}_{\Sigma}(X)}($};
\node[] () at (\posZerox,0) [] {$P_{0, 0}$};
\node[] (0c0) at (\posDotOnex,0) [] {$, \cdots,$};
\node[] () at (\posJx,0) [] {$P_{j, 0}$};
\node[] (0c1) at (\posDotTwox,0) [] {$, \cdots,$};
\node[] () at (\posFx,0) [] {$P_{\bb{\mathbf{s}}-1,0}$};
\node[] (1c0) at (\posEndx,0) [] {$)$};

\node[] (E0) at (\posPx,-.75) [] {\rotatebox{-90}{$\mor$}};
\node[] (E0j) at (\posZerox,-.75) [] {\rotatebox{-90}{$\mor$}};
\node[] (E00) at (\posZerox+.7,-.75) [] {$\scriptstyle (\mathfrak{p},T)_{0,0}$};
\node[] (E0j) at (\posJx,-.75) [] {\rotatebox{-90}{$=$}};
\node[] (E0f) at (\posFx,-.75) [] {\rotatebox{-90}{$\,=$}};

\node[] (1P) at (\posPx,-1.4) [] {$\vdots$};
\node[] () at (\posZerox,-1.4) [] {$\vdots$};
\node[] (1c0) at (\posDotOnex,-1.5) [] {$, \cdots,$};
\node[] (1j) at (\posJx,-1.4) [] {$\vdots$};
\node[] (1c1) at (\posDotTwox,-1.5) [] {$, \cdots,$};
\node[] (1f) at (\posFx,-1.4) [] {$\vdots$};

\node[] (E1) at (\posPx,-2.25) [] {\rotatebox{-90}{$\mor$}};
\node[] (E10) at (\posZerox,-2.25) [] {\rotatebox{-90}{$\mor$}};
\node[] (E10label) at (\posZerox+1,-2.25) [] {$\scriptstyle (\mathfrak{p},T)_{0, \bb{\mathbf{c}_{0}}-1}$};
\node[] (E1j) at (\posJx,-2.25) [] {\rotatebox{-90}{$\,=$}};
\node[] (E1f) at (\posFx,-2.25) [] {\rotatebox{-90}{$\,=$}};

\node[] (2P) at (\posPx,-3) [] {$P_{\bb{\mathbf{c}_{0}}}$};
\node[] (2=) at (\posEqx,-3) [] {$=$};
\node[] (2s) at (\posSigmax,-2.95) [] {$\sigma^{\mathbf{T}_{\Sigma}(X)}($};
\node[] () at (\posZerox,-3) [] {$P_{0,\bb{\mathbf{c}_{0}}}$};
\node[] (2c0) at (\posDotOnex,-3) [] {$, \cdots,$};
\node[] (2j) at (\posJx,-3) [] {$P_{j,0}$};
\node[] (2c1) at (\posDotTwox,-3) [] {$, \cdots,$};
\node[] (2f) at (\posFx,-3) [] {$P_{\bb{\mathbf{s}}-1, 0}$};
\node[] (2c0) at (\posEndx,-3) [] {$)$};

\node[] (E2) at (\posPx,-3.75) [] {\rotatebox{-90}{$\mor$}};
\node[] (E20) at (\posZerox,-3.75) [] {\rotatebox{-90}{$\,=$}};
\node[] (E2j) at (\posJx,-3.75) [] {\rotatebox{-90}{$\, =$}};
\node[] (E2f) at (\posFx,-3.75) [] {\rotatebox{-90}{$\,=$}};

\node[] (3P) at (\posPx,-4.4) [] {$\vdots$};
\node[] (30) at (\posZerox,-4.4) [] {$\vdots$};
\node[] (3c0) at (\posDotOnex,-4.5) [] {$\cdots$};
\node[] (3j) at (\posJx,-4.4) [] {$\vdots$};
\node[] (3c1) at (\posDotTwox,-4.5) [] {$\cdots$};
\node[] (3f) at (\posFx,-4.4) [] {$\vdots$};

\node[] (E3) at (\posPx,-5.25) [] {\rotatebox{-90}{$\mor$}};
\node[] (E30) at (\posZerox,-5.25) [] {\rotatebox{-90}{$\,=$}};
\node[] (E3j) at (\posJx,-5.25) [] {\rotatebox{-90}{$\,=$}};
\node[] (E3f) at (\posFx,-5.25) [] {\rotatebox{-90}{$\,=$}};

\node[] (4P) at (\posPx,-6) [] {$P_{\sum_{l=0}^{j-1}\bb{\mathbf{c}_{l}}}$};
\node[] (4=) at (\posEqx,-6) [] {$=$};
\node[] (4s) at (\posSigmax,-5.95) [] {$\sigma^{\mathbf{T}_{\Sigma}(X)}($};
\node[] () at (\posZerox,-6) [] {$P_{0, \bb{\mathbf{c}^{0}}}$};
\node[] (4c0) at (\posDotOnex,-6) [] {$, \cdots,$};
\node[] (4) at (\posJx,-6) [] {$P_{j,0}$};
\node[] (4c1) at (\posDotTwox,-6) [] {$, \cdots,$};
\node[] (4f) at (\posFx,-6) [] {$P_{\bb{\mathbf{s}}-1, 0}$};
\node[] (4c0) at (\posEndx,-6) [] {$)$};

\node[] (E4) at (\posPx,-6.75) [] {\rotatebox{-90}{$\mor$}};
\node[] (E40) at (\posZerox,-6.75) [] {\rotatebox{-90}{$\,=$}};
\node[] (E4j) at (\posJx,-6.75) [] {\rotatebox{-90}{$\mor$}};
\node[] (E4jlabel) at (\posJx+0.70, -6.75) [] {$\scriptstyle (\mathfrak{p},T)_{j, 0}$};
\node[] (E4f) at (\posFx,-6.75) [] {\rotatebox{-90}{$\,=$}};

\node[] (5P) at (\posPx,-7.4) [] {$\vdots$};
\node[] (50) at (\posZerox,-7.4) [] {$\vdots$};
\node[] (5c0) at (\posDotOnex,-7.5) [] {$\cdots$};
\node[] () at (\posJx,-7.4) [] {$\vdots$};
\node[] (5c1) at (\posDotTwox,-7.5) [] {$\cdots$};
\node[] (5f) at (\posFx,-7.5) [] {$\vdots$};

\node[] (E5) at (\posPx,-8.25) [] {\rotatebox{-90}{$\mor$}};
\node[] (E50) at (\posZerox,-8.25) [] {\rotatebox{-90}{$\,=$}};
\node[] () at (\posJx,-8.25) [] {\rotatebox{-90}{$\mor$}};
\node[] (E5jlabel) at (\posJx+1, -8.25) [] 
{$\scriptstyle (\mathfrak{p},T)_{j, \bb{\mathbf{c}_{j}}-1}$};
\node[] (E5f) at (\posFx,-8.25) [] {\rotatebox{-90}{$\,=$}};

\node[] (6P) at (\posPx,-9) [] {$P_{\sum^{j}_{l=0}\bb{\mathbf{c}_{l}}}$};
\node[] (6=) at (\posEqx,-9) [] {$=$};
\node[] (6s) at (\posSigmax,-8.95) [] {$\sigma^{\mathbf{T}_{\Sigma}(X)}($};
\node[] () at (\posZerox,-9) []
{$\color{white}P_{\bb{\mathbf{c}_{\bb{\mathbf{s}}-1}}, \bb{\mathbf{s}}-1}$};
\node[] () at (\posZerox,-9) [] {$P_{0, \bb{\mathbf{c}_{0}}}$};
\node[] (6c0) at (\posDotOnex,-9) [] {$, \cdots,$};
\node[] () at (\posJx,-9) [] {$P_{j, \bb{\mathbf{c}_{j}}}$};
\node[] (6c1) at (\posDotTwox,-9) [] {$, \cdots,$};
\node[] (6f) at (\posFx,-9) [] {$P_{\bb{\mathbf{s}}-1, 0}$};
\node[] (6c0) at (\posEndx,-9) [] {$)$};

\node[] (E6) at (\posPx,-9.75) [] {\rotatebox{-90}{$\mor$}};
\node[] (E60) at (\posZerox,-9.75) [] {\rotatebox{-90}{$\,=$}};
\node[] (E6j) at (\posJx,-9.75) [] {\rotatebox{-90}{$\,=$}};
\node[] (E6f) at (\posFx,-9.75) [] {\rotatebox{-90}{$\,=$}};

\node[] (7P) at (\posPx,-10.25) [] {$\vdots$};
\node[] (70) at (\posZerox,-10.25) [] {$\vdots$};
\node[] (7c0) at (\posDotOnex,-10.25) [] {$\cdots$};
\node[] (7j) at (\posJx,-10.25) [] {$\vdots$};
\node[] (7c1) at (\posDotTwox,-10.25) [] {$\cdots$};
\node[] (7f) at (\posFx,-10.25) [] {$\vdots$};

\node[] (E7) at (\posPx,-11) [] {\rotatebox{-90}{$\mor$}};
\node[] (E70) at (\posZerox,-11) [] {\rotatebox{-90}{$\,=$}};
\node[] (E7j) at (\posJx,-11) [] {\rotatebox{-90}{$\,=$}};
\node[] (E7f) at (\posFx,-11) [] {\rotatebox{-90}{$\,=$}};

\node[] (8P) at (\posPx,-11.75) [] {$P_{\sum_{l=0}^{\bb{\mathbf{s}}-2}\bb{\mathbf{c}_{l}}}$};
\node[] (8=) at (\posEqx,-11.75) [] {$=$};
\node[] (8s) at (\posSigmax,-11.7) [] {$\sigma^{\mathbf{T}_{\Sigma}(X)}($};
\node[] () at (\posZerox,-11.75) [] {$P_{0, \bb{\mathbf{c}_{0}}}$};
\node[] (8c0) at (\posDotOnex,-11.75) [] {$, \cdots,$};
\node[] () at (\posJx,-11.75) [] {$P_{j, \bb{\mathbf{c}_{j}}}$};
\node[] (8c1) at (\posDotTwox,-11.75) [] {$, \cdots,$};
\node[] (8f) at (\posFx,-11.75) [] {$P_{\bb{\mathbf{s}}-1, 0}$};
\node[] (8c0) at (\posEndx,-11.75) [] {$)$};

\node[] (E8) at (\posPx,-12.5) [] {\rotatebox{-90}{$\mor$}};
\node[] (E80) at (\posZerox,-12.5) [] {\rotatebox{-90}{$\,=$}};
\node[] (E8j) at (\posJx,-12.5) [] {\rotatebox{-90}{$\,=$}};
\node[] (E8flabel) at (\posFx+.95, -12.5) [] 
{$\scriptstyle (\mathfrak{p},T)_{\bb{\mathbf{s}}-1, 0}$};
\node[] (E8f) at (\posFx,-12.5) [] {\rotatebox{-90}{$\mor$}};

\node[] (9P) at (\posPx,-13.24) [] {$\vdots$};
\node[] (90) at (\posZerox,-13.24) [] {$\vdots$};
\node[] (9c0) at (\posDotOnex,-13.25) [] {$\cdots$};
\node[] (9j) at (\posJx,-13.24) [] {$\vdots$};
\node[] (9c1) at (\posDotTwox,-13.25) [] {$\cdots$};
\node[] () at (\posFx,-13.24) [] {$\vdots$};

\node[] (E9) at (\posPx,-14) [] {\rotatebox{-90}{$\mor$}};
\node[] (E90) at (\posZerox,-14) [] {\rotatebox{-90}{$\,=$}};
\node[] (E9j) at (\posJx,-14) [] {\rotatebox{-90}{$\,=$}};
\node[] (E9flabel) at (\posFx+1.45, -14) [] 
{$\scriptstyle (\mathfrak{p},T)_{\bb{\mathbf{s}}-1, \bb{\mathbf{c}_{\bb{\mathbf{s}}-1}}-1}$};
\node[] (E9f) at (\posFx,-14) [] {\rotatebox{-90}{$\mor$}};

\node[] (10P) at (\posPx,-14.75) [] {$P_{\bb{\mathbf{c}}}$};
\node[] (10=) at (\posEqx,-14.75) [] {$=$};
\node[] (10s) at (\posSigmax,-14.7) [] {$\sigma^{\mathbf{T}_{\Sigma}(X)}($};
\node[] () at (\posZerox,-14.75) [] {$P_{0, \bb{\mathbf{c}_{0}}}$};
\node[] (10c0) at (\posDotOnex,-14.75) [] {$, \cdots,$};
\node[] () at (\posJx,-14.75) [] {$P_{j, \bb{\mathbf{c}_{j}}}$};
\node[] (10c1) at (\posDotTwox,-14.75) [] {$, \cdots,$};
\node[] (10f) at (\posFx,-14.75) [] {$P_{\bb{\mathbf{s}}-1, \bb{\mathbf{c}_{\bb{\mathbf{s}}-1}}}$};
\node[] (10c0) at (\posEndx,-14.75) [] {$)$};
\end{tikzpicture}
\end{scaletikzpicturetowidth}
\caption{The path $\sigma^{\mathbf{Pth}_{\boldsymbol{\mathcal{A}}}}((\mathfrak{P}_{j})_{j\in\bb{\mathbf{s}}})$.}\label{FPthSigma}
\end{figure}

In the following corollary we state that a one-step echelonless path can always be represented as the action of an univocally determined operation on the family of paths extracted from it.

\begin{restatable}{corollary}{CUStep}
\label{CUStep}
Let $s$ be a sort in $S$, $\mathfrak{P}$ a one-step echelonless path in $\mathrm{Pth}_{\boldsymbol{\mathcal{A}},s}$ of type $\sigma\in \Sigma_{\mathbf{s},s}$, for a unique pair $(\mathbf{s},s)\in (S^{\star}-\{\lambda\})\times S$, and $(\mathfrak{P}_{j})_{j\in \bb{\mathbf{s}}}$ the family of  paths we can extract from $\mathfrak{P}$ in virtue of Lemma~\ref{LPthExtract}. Then we have that
\[
\mathfrak{P}=\sigma^{\mathbf{Pth}_{\boldsymbol{\mathcal{A}}}}\left(\left(
\mathfrak{P}_{j}\right)_{j\in\bb{\mathbf{s}}}\right).
\] 
\end{restatable}

\begin{proof}
Let us assume that, for a sort $t\in S$, $\mathfrak{P}$ is a  $t$-path in $\mathrm{Pth}_{\boldsymbol{\mathcal{A}},s}$ of the form
\[
\mathfrak{P}=
\left(
\left(
P_{i}
\right)_{i\in 2},
\mathfrak{p},
T
\right),
\]
where $T$ is a $t$-translation of sort $s$. Then, by Definition~\ref{DTrans}, for an $n\in\mathbb{N}-\{0\}$, $T$ is a composition of elementary translations $T_{n-1}\circ \cdots\circ  T_{0}$. Let $T'$ be the maximal prefix of $T$.

If $T$ is a translation associated to the operation symbol $\sigma\in\Sigma_{\mathbf{s},s}$, then there exists an index $k\in\bb{\mathbf{s}}$, a family of terms $(P_{j})_{j\in k}\in \prod_{j\in k}\mathrm{T}_{\Sigma}(X)_{s_{j}}$ and a family of terms $(P_{l})_{l\in \bb{\mathbf{s}}-(k+1)}\in \prod_{l\in \bb{\mathbf{s}}-(k+1)}\mathrm{T}_{\Sigma}(X)_{s_{l}}$ such that $s_{k}=t$ and
\[
T=\sigma^{\mathbf{T}_{\Sigma^{\mathrm{c},\mathcal{A}}}(X)}\left(
P_{0},
\cdots,
P_{k-1},
T',
P_{k+1},
\cdots,
P_{\bb{\mathbf{s}}-1}
\right).
\]

If $\mathfrak{p}=(M,N)$, then, since $\mathfrak{P}$ is a path, we have that 
\begin{multicols}{2}
\begin{itemize}
\item[(i)] $T(M)=P_{0}$;
\item[(ii)]  $T(N)=P_{1}$.
\end{itemize}
\end{multicols}

By Lemma~\ref{LPthExtract}, if $(\mathfrak{P})_{j\in\bb{\mathbf{s}}}$ is the family of paths we can extract from $\mathfrak{P}$, then we have that 
\[
\mathfrak{P}_{j}=
\begin{cases}
\mathrm{ip}^{(1,0)\sharp}_{s_{j}}\left(
P_{j}
\right)
&\qquad\mbox{if } j\neq k;
\\
\mathfrak{P}_{k}
&\qquad\mbox{if } j=k,
\end{cases}
\] 
where, by Lemma~\ref{LPthExtract}, we have that
$$
\xymatrix@C=50pt{
\mathfrak{P}_{k}\colon
T'(M)
\ar@{->}[r]^-{\text{\Small{($\mathfrak{p}$, $
T'
$)}}}
&
T'(N)
}
.
$$

Now, if we consider the interpretation of the operation symbol $\sigma$ in the $S$-sorted $\Sigma$-algebra $\mathbf{Pth}_{\boldsymbol{\mathcal{A}}}$, introduced in Proposition~\ref{PPthAlg}, then we have that 
$\sigma^{\mathbf{Pth}_{\boldsymbol{\mathcal{A}}}}((\mathfrak{P}_{j})_{j\in\bb{\mathbf{s}}})$ is a  $t$-path in 
$\mathrm{Pth}_{\boldsymbol{\mathcal{A}},s}$ of the form
\[
\sigma^{\mathbf{Pth}_{\boldsymbol{\mathcal{A}}}}\left(\left(\mathfrak{P}_{j}
\right)_{j\in\bb{\mathbf{s}}}\right)
=
\left(
\left(
Q_{i}
\right)_{i\in 2}
,
\mathfrak{p}
,
U
\right),
\]

Regarding the initial $0$-constituent of $\sigma^{\mathbf{Pth}_{\boldsymbol{\mathcal{A}}}}((\mathfrak{P}_{j})_{j\in\bb{\mathbf{s}}})$,  the following chain of equalities holds
\allowdisplaybreaks
\begin{align*}
Q_{0}&=
\sigma^{\mathbf{T}_{\Sigma}(X)}
\left(
P_{0}
,
\cdots,
P_{k-1},
T'\left(
M
\right),
P_{k+1}
\cdots,
P_{\bb{\mathbf{s}}-1}
\right)
\tag{1}
\\&=
T\left(
M
\right)
\tag{2}
\\&=
P_{0}.
\tag{3}
\end{align*}

In the just stated chain of equalities, the first equality unravels the description of the initial $0$-constituent of $\sigma^{\mathbf{Pth}_{\boldsymbol{\mathcal{A}}}}((\mathfrak{P}_{j})_{j\in\bb{\mathbf{s}}})$, by  Proposition~\ref{PPthAlg}; the second equality recovers the description of the elementary  translation $T$; finally, the last equality follows from item (i) above.

Regarding the final $0$-constituent of $\sigma^{\mathbf{Pth}_{\boldsymbol{\mathcal{A}}}}((\mathfrak{P}_{j})_{j\in\bb{\mathbf{s}}})$,  the following chain of equalities holds
\allowdisplaybreaks
\begin{align*}
Q_{1}&=
\sigma^{\mathbf{T}_{\Sigma}(X)}
\left(
P_{0},
\cdots,
P_{k-1},
T'\left(
N
\right),
P_{k+1},
\cdots,
P_{\bb{\mathbf{s}}-1}
\right)
\tag{1}
\\&=
T\left(
N
\right)
\tag{2}
\\&=
P_{1}.
\tag{3}
\end{align*}

In the just stated chain of equalities, the first equality unravels the description of the final $0$-constituent of $\sigma^{\mathbf{Pth}_{\boldsymbol{\mathcal{A}}}}((\mathfrak{P}_{j})_{j\in\bb{\mathbf{s}}})$, by  Proposition~\ref{PPthAlg}; the second equality recovers the description of the elementary  translation $T$; finally, the last equality follows from item (ii) above.

Finally,  $U$,  the unique translation occurring in  $\sigma^{\mathbf{Pth}_{\boldsymbol{\mathcal{A}}}}((\mathfrak{P}_{j})_{j\in\bb{\mathbf{s}}})$, satisfies the following chain of equalities 
\allowdisplaybreaks
\begin{align*}
U&=
\sigma^{\mathbf{T}_{\Sigma}(X)}\left(
P_{0},
\cdots,
P_{k-1},
T',
P_{k+1},
\cdots,
P_{\bb{\mathbf{s}}-1}
\right)
\tag{1}
\\&=
T.\tag{2}
\end{align*}

In the just stated chain of equalities, the first equality unravels the description of the unique  translation occurring in  $\sigma^{\mathbf{Pth}_{\boldsymbol{\mathcal{A}}}}((\mathfrak{P}_{j})_{j\in\bb{\mathbf{s}}})$, by Proposition~\ref{PPthAlg}; finally, the last equality recovers the description of the translation $T$.

All in all, we conclude that 
\[
\mathfrak{P}=\sigma^{\mathbf{Pth}_{\boldsymbol{\mathcal{A}}}}\left(\left(
\mathfrak{P}_{j}\right)_{j\in\bb{\mathbf{s}}}\right).
\]

This completes the proof.
\end{proof}

\section{
\texorpdfstring
{On the subalgebra of identity paths}
{The subalgebra of identity paths}
}

We next study how the previous operations behave when acting on $(1,0)$-identity paths. We first note that the interpretation of every constant in $\Sigma$ is a $(1,0)$-identity path.

\begin{restatable}{proposition}{RConsSigma}
\label{RConsSigma}
Let $s$ be a sort in $S$ and $\sigma$ a constant operation symbol in $\Sigma_{\lambda,s}$. Then its realization as a constant in $\mathbf{Pth}^{(0,1)}_{\boldsymbol{\mathcal{A}}}$, i.e., $\sigma^{\mathbf{Pth}_{\boldsymbol{\mathcal{A}}}}$, is $(\sigma^{\mathbf{T}_{\Sigma}(X)},\lambda,\lambda)$, the $(1,0)$-identity path on the term $\sigma^{\mathbf{T}_{\Sigma}(X)}$.
\end{restatable}

We next prove that every operation of the $S$-sorted $\Sigma$-algebra $\mathbf{Pth}_{\boldsymbol{\mathcal{A}}}^{(0,1)}$ when restricted to $(1,0)$-identity paths retrieves $(1,0)$-identity paths.

\begin{restatable}{proposition}{PSigma}
\label{PSigma} Let $(\mathbf{s},s)$ be a pair in $(S^{\star}-\{\lambda\})\times S$, $\sigma$ an operation symbol in $\Sigma_{\mathbf{s},s}$ and $(\mathfrak{P}_{j})$ a family of $(1,0)$-identity paths in 
$\mathrm{Pth}_{\boldsymbol{\mathcal{A}},s}$, where, for every $j\in\bb{\mathbf{s}}$, $\mathfrak{P}_{j}=\mathrm{ip}^{(1,0)\sharp}_{s_{j}}(P_{j})$, for a suitable family of terms $(P_{j})_{j\in\bb{\mathbf{s}}}$ in 
$\mathrm{T}_{\Sigma}(X)_{\mathbf{s}}$. Then
$$
\sigma^{\mathbf{Pth}_{\boldsymbol{\mathcal{A}}}}
\left(\left(
\mathfrak{P}_{j}
\right)_{j\in\bb{\mathbf{s}}}
\right)
=
\mathrm{ip}^{(1,0)\sharp}_{s}
\left(
\sigma^{\mathbf{T}_{\Sigma}(X)}
\left(
\left(
P_{j}
\right)_{j\in\bb{\mathbf{s}}}
\right)
\right).
$$ 
\end{restatable}

\begin{proof}
From Claim~\ref{CPthSigma} it follows that $\sigma^{\mathbf{Pth}_{\boldsymbol{\mathcal{A}}}}
((
\mathfrak{P}_{j}
)_{j\in\bb{\mathbf{s}}}
)$ is a path of length $0$ on the term $\sigma^{\mathbf{T}_{\Sigma}(X)}
((
P_{j}
)_{j\in\bb{\mathbf{s}}}
)$.
\end{proof}

\begin{restatable}{corollary}{RSigma}\label{RSigma}
The $S$-sorted subset $\mathrm{ip}^{(1,0)\sharp}[\mathrm{T}_{\Sigma}(X)]$, of the $(1,0)$-identity paths of $\mathrm{Pth}_{\boldsymbol{\mathcal{A}}}$, is a closed subset of $\mathbf{Pth}_{\boldsymbol{\mathcal{A}}}^{(0,1)}$.
\end{restatable}

\begin{restatable}{definition}{DIpz}
\label{DIpz}\index{identity!first-order!$\mathrm{ip}^{(1,0)\sharp}[\mathbf{T}_{\Sigma}(X)]$}
We let $\mathrm{ip}^{(1,0)\sharp}[\mathbf{T}_{\Sigma}(X)]$ stand for the $\Sigma$-subalgebra of $\mathbf{Pth}^{(0,1)}_{\boldsymbol{\mathcal{A}}}$ whose underlying $S$-sorted set is $\mathrm{ip}^{(1,0)\sharp}[\mathrm{T}_{\Sigma}(X)]$.
\end{restatable}

\section{
\texorpdfstring
{Connecting layers $0$ and $1$}
{Connecting layers 0 and 1}
}
In this section, we investigate how the many-sorted algebras of layers $1$ and $0$ are interconnected.


We begin by proving that the mappings $\mathrm{sc}^{(0,1)}$ and $\mathrm{tg}^{(0,1)}$ are $\Sigma$-homomorphisms from $\mathbf{Pth}_{\boldsymbol{\mathcal{A}}}^{(0,1)}$ to $\mathbf{T}_{\Sigma}(X)$.

\begin{restatable}{proposition}{PHom}
\label{PHom} The mappings $\mathrm{sc}^{(0,1)}$ and $\mathrm{tg}^{(0,1)}$ are $\Sigma$-homomorphisms from $\mathbf{Pth}^{(0,1)}_{\boldsymbol{\mathcal{A}}}$ to $\mathbf{T}_{\Sigma}(X)$.
\end{restatable}
\begin{proof}
Let $(\mathbf{s},s)$ be an element of $S^{\star}\times S$, $\sigma$ an operation symbol in $\Sigma_{\mathbf{s},s}$ and $(\mathfrak{P}_{j})_{j\in\bb{\mathbf{s}}}$ a family of paths in $\mathrm{Pth}_{\boldsymbol{\mathcal{A}},\mathbf{s}}$, then according to Claim~\ref{CPthSigma}, we have that
\begin{align*}
\mathrm{sc}^{(0,1)}_{s}\left(
\sigma^{\mathbf{Pth}_{\boldsymbol{\mathcal{A}}}}
\left(\left(\mathfrak{P}_{j}
\right)_{j\in\bb{\mathbf{s}}}
\right)\right)
&=
\sigma^{\mathbf{T}_{\Sigma}(X)}
\left(\left(
\mathrm{sc}^{(0,1)}_{s_{j}}
\left(\mathfrak{P}_{j}
\right)
\right)_{j\in\bb{\mathbf{s}}}\right),
\\
\mathrm{tg}^{(0,1)}_{s}\left(
\sigma^{\mathbf{Pth}_{\boldsymbol{\mathcal{A}}}}
\left(\left(\mathfrak{P}_{j}
\right)_{j\in\bb{\mathbf{s}}}
\right)
\right)
&=
\sigma^{\mathbf{T}_{\Sigma}(X)}
\left(\left(
\mathrm{tg}^{(0,1)}_{s_{j}}
\left(\mathfrak{P}_{j}
\right)
\right)_{j\in\bb{\mathbf{s}}}
\right).
\end{align*}
The above statement also includes the particular case of constant symbols in $\Sigma_{\lambda,s}$ according to Remark~\ref{RConsSigma}.

This proves Proposition~\ref{PHom}.
\end{proof}

Also $\mathrm{ip}^{(1,0)\sharp}$, the $(1,0)$-identity path mapping, is a $\Sigma$-homomorphism from $\mathbf{T}_{\Sigma}(X)$ to $\mathbf{Pth}_{\boldsymbol{\mathcal{A}}}^{(0,1)}$. To obtain this result we will make use of the universal property of the free $S$-sorted $\Sigma$-algebra $\mathbf{T}_{\Sigma}(X)$.

\begin{restatable}{proposition}{PIpHom}
\label{PIpHom}
The mapping $\mathrm{ip}^{(1,0)\sharp}$ is a $\Sigma$-homomorphism from $\mathbf{T}_{\Sigma}(X)$ to $\mathbf{Pth}^{(0,1)}_{\boldsymbol{\mathcal{A}}}$.
\end{restatable}
\begin{proof}
Let us recall from Definition~\ref{DPth} that $\mathrm{ip}^{(1,X)}$ is the $S$-sorted mapping from $X$ to $\mathbf{Pth}^{(0,1)}_{\boldsymbol{\mathcal{A}}}$ that, for every sort $s\in S$, sends $x\in X_{s}$ to $(x,\lambda,\lambda)$, the $(1,0)$-identity path on $x$.

By Proposition~\ref{PPthAlg}, we have that $\mathbf{Pth}^{(0,1)}_{\boldsymbol{\mathcal{A}}}$ is a $\Sigma$-algebra. By the universal property of the free $\Sigma$-algebra $\mathbf{T}_{\Sigma}(X)$, there exists a unique $\Sigma$-homomorphism $\mathrm{ip}^{(1,0)\sharp}$ from $\mathbf{T}_{\Sigma}(X)$ to $\mathbf{Pth}^{(0,1)}_{\boldsymbol{\mathcal{A}}}$ such that
$$
\mathrm{ip}^{(1,0)\sharp}\circ\eta^{(0,X)}=\mathrm{ip}^{(1,X)}
$$
i.e., such that the diagram in Figure~\ref{FHom} commutes.

Let us recall that the $\Sigma$-homomorphism $\mathrm{ip}^{(1,0)\sharp}$ sends, for every sort $s\in S$, a term $P\in\mathrm{T}_{\Sigma}(X)_{s}$ to $(P,\lambda,\lambda)$, the $(1,0)$-identity path on $P$.
\end{proof}

\begin{figure}
\begin{tikzpicture}
[ACliment/.style={-{To [angle'=45, length=5.75pt, width=4pt, round]}},scale=1]
\node[] (xoq) at (0,0) [] {$X$};
\node[] (txoq) at (6,0) [] {$\mathbf{T}_{\Sigma}(X)$};
\node[] (p) at (6,-3) [] {$\mathbf{Pth}^{(0,1)}_{\boldsymbol{\mathcal{A}}}$};
\draw[ACliment]  (xoq) to node [above]
{$\eta^{(0,X)}$} (txoq);
\draw[ACliment, bend right=10]   (xoq) to node [below left] {$\mathrm{ip}^{(1,X)}$} (p);

\node[] (B0) at (6,-1.5)  [] {};
\draw[ACliment]  ($(B0)+(0,1.2)$) to node [above, fill=white] {
$\textstyle \mathrm{ip}^{(1,0)\sharp}$
} ($(B0)+(0,-1.2)$);
\draw[ACliment, bend right]  ($(B0)+(.3,-1.2)$) to node [ below, fill=white] {
$\textstyle \mathrm{tg}^{(0,1)}$
} ($(B0)+(.3,1.2)$);
\draw[ACliment, bend left]  ($(B0)+(-.3,-1.2)$) to node [below, fill=white] {
$\textstyle \mathrm{sc}^{(0,1)}$
} ($(B0)+(-.3,1.2)$);
\end{tikzpicture}
\caption{$\Sigma$-homomorphisms relative to $X$ at layers $0$ \& $1$.}\label{FHom}
\end{figure}

From the above results, we conclude that the free $\Sigma$-algebra of terms is isomorphic to the many-sorted $\Sigma$-subalgebra of $(1,0)$-identity paths.

\begin{corollary}\label{CIsoIp} The $\Sigma$-homomorphism $\mathrm{ip}^{(1,0)\sharp}$ is a $\Sigma$-isomorphism from $\mathbf{T}_{\Sigma}(X)$ to $\mathrm{ip}^{(1,0)\sharp}[\mathbf{T}_{\Sigma}(X)]$.
\end{corollary}
\begin{proof}
By Proposition~\ref{PIpHom}, $\mathrm{ip}^{(1,0)\sharp}$ is a $\Sigma$-homomorphism from $\mathbf{T}_{\Sigma}(X)$ to $\mathbf{Pth}_{\boldsymbol{\mathcal{A}}}^{(0,1)}$ that correstricts to $\mathrm{ip}^{(1,0)\sharp}[\mathbf{T}_{\Sigma}(X)]$. That $\mathrm{ip}^{(1,0)\sharp}$ is surjective follows from the description of the codomain and injective by Proposition~\ref{PBasicEq}, since $\mathrm{sc}^{(0,1)}\circ\mathrm{ip}^{(1,0)\sharp}=\mathrm{id}^{\mathbf{T}_{\Sigma}(X)}$ (analogously with the $(0,1)$-target).
\end{proof}

\begin{proposition}\label{PSection}
The following equations between $\Sigma$-homomorphisms hold
\begin{multicols}{2}
\begin{itemize}
\item[(i)] $\mathrm{sc}^{(0,1)}\circ \mathrm{ip}^{(1,0)\sharp}
=\mathrm{id}^{\mathbf{T}_{\Sigma}(X)}$,
\item[(ii)] $\mathrm{tg}^{(0,1)}\circ \mathrm{ip}^{(1,0)\sharp}
=\mathrm{id}^{\mathbf{T}_{\Sigma}(X)}$.
\end{itemize}
\end{multicols}
\end{proposition}
\begin{proof}
This follows from Proposition~\ref{PBasicEq}.
\end{proof}

Thus, $\mathbf{T}_{\Sigma}(X)$ is a retract of $\mathbf{Pth}^{(0,1)}_{\boldsymbol{\mathcal{A}}}$, the
$\Sigma$-homomorphisms $\mathrm{sc}^{(0,1)}$ and $\mathrm{tg}^{(0,1)}$ are retractions of $\mathbf{Pth}^{(0,1)}_{\boldsymbol{\mathcal{A}}}$ onto $\mathbf{T}_{\Sigma}(X)$, and the $\Sigma$-homomorphism $\mathrm{ip}^{(1,0)\sharp}$ is a section of both $\mathrm{sc}^{(0,1)}$ and $\mathrm{tg}^{(0,1)}$.

\begin{corollary}\label{CSgTg}
Let $s$ be a sort in $S$ and $\mathfrak{P}$ a path in $\mathrm{Pth}_{\boldsymbol{\mathcal{A}},s}$. Then we have that
\begin{align*}
\left(\mathfrak{P},\mathrm{ip}^{(1,0)\sharp}_{s}\left(\mathrm{sc}^{(0,1)}_{s}\left(\mathfrak{P}\right)\right)\right)
&\in\mathrm{Ker}\left(\mathrm{sc}^{(0,1)}\right)_{s}\text{ and}\\
\left(\mathfrak{P},\mathrm{ip}^{(1,0)\sharp}_{s}\left(\mathrm{tg}^{(0,1)}_{s}\left(\mathfrak{P}\right)\right)\right)
&\in\mathrm{Ker}\left(\mathrm{tg}^{(0,1)}\right)_{s}.
\end{align*}
\end{corollary}

\begin{remark}
From Proposition~\ref{PHom} it follows that the kernel of $\mathrm{sc}^{(0,1)}$ is a $\Sigma$-congruence on $\mathbf{Pth}^{(0,1)}_{\boldsymbol{\mathcal{A}}}$. Therefore, we can consider the standard (epi, mono)-factoriza\-tion of $\mathrm{sc}^{(0,1)}$ given by
\begin{enumerate}
\item the projection $\mathrm{pr}^{\mathrm{Ker}(\mathrm{sc}^{(0,1)})}$ from $\mathbf{Pth}^{(0,1)}_{\boldsymbol{\mathcal{A}}}$ to $\mathbf{Pth}^{(0,1)}_{\boldsymbol{\mathcal{A}}}/{\mathrm{Ker}(\mathrm{sc^{(0,1)}})}$ that, for every $s\in S$, assigns to a path $\mathfrak{P}$ in $\mathrm{Pth}_{\boldsymbol{\mathcal{A}},s}$ its equivalence class $[\mathfrak{P}]_{\mathrm{Ker}(\mathrm{sc}^{(0,1)})_{s}}$; and
\item the embedding $\mathrm{sc}^{(0,1)\mathrm{m}}$ of $\mathbf{Pth}^{(0,1)}_{\boldsymbol{\mathcal{A}}}/{\mathrm{Ker}(\mathrm{sc^{(0,1)}})}$ into $\mathbf{T}_{\Sigma}(X)$ that, for every $s\in S$, assigns to an equivalence class $[\mathfrak{P}]_{\mathrm{Ker}(\mathrm{sc}^{(0,1)})_{s}}$ in $\mathrm{Pth}_{\boldsymbol{\mathcal{A}},s}/{\mathrm{Ker}(\mathrm{sc}^{(0,1)})_{s}}$, with $\mathfrak{P}\in\mathrm{Pth}_{\boldsymbol{\mathcal{A}},s}$, the term $\mathrm{sc}^{(0,1)}_{s}(\mathfrak{P})$, i.e., the value of the mapping $\mathrm{sc}^{(0,1)}_{s}$ at any equivalence class representative.  
We will refer to the mapping $\mathrm{sc}^{(0,1)\mathrm{m}}$ as the \emph{monomorphic factorization of $\mathrm{sc}^{(0,1)}$}.
\end{enumerate}  
Let us recall that $(\mathrm{pr}^{\mathrm{Ker}(\mathrm{sc}^{(0,1)})}, \mathrm{sc}^{(0,1)\mathrm{m}})$ is the unique pair of $\Sigma$-homomorphisms such that $\mathrm{sc}^{(0,1)\mathrm{m}}\circ\mathrm{pr}^{\mathrm{Ker}(\mathrm{sc^{(0,1)}})}=\mathrm{sc}^{(0,1)}$, i.e., such that the upper half of the diagram in Figure~\ref{FIsos} commutes.

On the other hand, the kernel of the $\mathrm{tg}^{(0,1)}$ is also a $\Sigma$-congruence on $\mathbf{Pth}^{(0,1)}_{\boldsymbol{\mathcal{A}}}$. We can consider the standard (epi, mono)-factorization of $\mathrm{tg}^{(0,1)}$ given by
\begin{enumerate}
\item the projection $\mathrm{pr}^{\mathrm{Ker}(\mathrm{tg}^{(0,1)})}$ from $\mathbf{Pth}^{(0,1)}_{\boldsymbol{\mathcal{A}}}$ to $\mathbf{Pth}^{(0,1)}_{\boldsymbol{\mathcal{A}}}/{\mathrm{Ker}(\mathrm{tg}^{(0,1)})}$ that, for every $s\in S$, assigns to a path $\mathfrak{P}$ in $\mathrm{Pth}_{\boldsymbol{\mathcal{A}},s}$ its equivalence class $[\mathfrak{P}]_{\mathrm{Ker}(\mathrm{tg}^{(0,1)})_{s}}$; and
\item the embedding $\mathrm{tg}^{(0,1)\mathrm{m}}$ from $\mathbf{Pth}^{(0,1)}_{\boldsymbol{\mathcal{A}}}/{\mathrm{Ker}(\mathrm{tg}^{(0,1)})}$ to $\mathbf{T}_{\Sigma}(X)$ that, for every $s\in S$, assigns to an equivalence class $[\mathfrak{P}]_{\mathrm{Ker}(\mathrm{tg}^{(0,1)})_{s}}$ in $\mathrm{Pth}_{\boldsymbol{\mathcal{A}},s}/{\mathrm{Ker}(\mathrm{tg}^{(0,1)})_{s}}$ and $\mathfrak{P}\in\mathrm{Pth}_{\boldsymbol{\mathcal{A}},s}$ the term $\mathrm{tg}^{(0,1)}_{s}(\mathfrak{P})$, i.e., the value of the mapping $\mathrm{tg}^{(0,1)}_{s}$ at any equivalence class representative. We will refer to the mapping $\mathrm{tg}^{(0,1)\mathrm{m}}$ as the \emph{monomorphic factorization of $\mathrm{tg}^{(0,1)}$}.
\end{enumerate}
Let us recall that $(\mathrm{pr}^{\mathrm{Ker}(\mathrm{tg}^{(0,1)})}, \mathrm{tg}^{(0,1)\mathrm{m}})$ is the unique pair of $\Sigma$-homomorphisms such that $\mathrm{tg}^{(0,1)\mathrm{m}}\circ\mathrm{pr}^{\mathrm{Ker}(\mathrm{tg}^{(0,1)})}=\mathrm{tg}^{(0,1)}$, i.e., such that the lower half of diagram in Figure~\ref{FIsos} commutes.
\end{remark}

\begin{figure}
\begin{tikzpicture}
[ACliment/.style={-{To [angle'=45, length=5.75pt, width=4pt, round]}}]
\node[] (txoq) at (-1,0) [] {$\mathbf{T}_{\Sigma}(X)$};
\node[] (p) at (4,3.5) [] {$\mathbf{Pth}^{(0,1)}_{\boldsymbol{\mathcal{A}}}$};
\node[] (p2) at (4,-3.5) [] {$\mathbf{Pth}^{(0,1)}_{\boldsymbol{\mathcal{A}}}$};
\node[] (pkernel) at (4,1.5) [] {$\mathbf{Pth}^{(0,1)}_{\boldsymbol{\mathcal{A}}}/{\mathrm{Ker}(\mathrm{sc}^{(0,1)})}$};
\node[] (pkernel2) at (4,-1.5) [] {$\mathbf{Pth}^{(0,1)}_{\boldsymbol{\mathcal{A}}}/{\mathrm{Ker}(\mathrm{tg}^{(0,1)})}$};
\node[] (txoq2) at (9,0) [] {$\mathbf{T}_{\Sigma}(X)$};
\draw[ACliment, bend left=15]  (txoq) to node [above left]
{$\mathrm{ip}^{(1,0)\sharp}$} (p);
\draw[ACliment]  (txoq) to node [sloped, above] {} (pkernel);
\draw[ACliment]  (p) to node [midway, fill=white]
{$\mathrm{pr}^{\mathrm{Ker}(\mathrm{sc}^{(0,1)})}$} (pkernel);
\draw[ACliment]  (pkernel) to node [sloped, above]
{$\mathrm{sc}^{(0,1)\mathrm{m}}$} (txoq2);
\draw[ACliment, bend left=15]  (p) to node [above right]
{$\mathrm{sc}^{(0,1)}$} (txoq2);
\draw[ACliment]  (txoq) to node  [midway, fill=white]
{$\mathrm{id}^{\mathbf{T}_{\Sigma}(X)}$} (txoq2);
\draw[ACliment, bend right=15]  (txoq) to node [below left]
{$\mathrm{ip}^{(1,0)\sharp}$} (p2);
\draw[ACliment]  (txoq) to node [sloped, below] {} (pkernel2);
\draw[ACliment]  (p2) to node  [midway, fill=white]
{$\mathrm{pr}^{\mathrm{Ker}(\mathrm{tg}^{(0,1)})}$} (pkernel2);
\draw[ACliment]  (pkernel2) to node [sloped, below]
{$\mathrm{tg}^{(0,1)\mathrm{m}}$} (txoq2);
\draw[ACliment, bend right=15]  (p2) to node [below right]
{$\mathrm{tg}^{(0,1)}$} (txoq2);
\end{tikzpicture}
\caption{$\Sigma$-algebras and $\Sigma$-homomorphisms at layers $1$ \& $0$.}
\label{FIsos}
\end{figure}

By the first isomorphism theorem, the aforementioned $\Sigma$-algebras are isomorphic.

\begin{restatable}{corollary}{CIsos}
\label{CIsos}
The $S$-sorted $\Sigma$-algebras 
\begin{multicols}{3}
\begin{itemize}
\item[(i)] $\mathbf{Pth}^{(0,1)}_{\boldsymbol{\mathcal{A}}}/{\mathrm{Ker}(\mathrm{sc}^{(0,1)})}$,
\item[(ii)] $\mathbf{Pth}^{(0,1)}_{\boldsymbol{\mathcal{A}}}/{\mathrm{Ker}(\mathrm{tg}^{(0,1)})}$,
\item[(iii)] $\mathbf{T}_{\Sigma}(X)$.
\end{itemize}
\end{multicols}
are isomorphic.
\end{restatable}

For the sake of illustration, here is an example of the notions defined in this chapter. From this example, and as if it were a \emph{roman-fleuve}, we will provide, in some of the subsequent chapters, further examples that will serve to illuminate the new concepts contained in them.

\begin{example}\label{ERun} 
Let us take $2=\{0,1\}$ as our set of sorts. Then we denote by $\Sigma$ the $2$-sorted signature, defined, for every $(\mathbf{s},s)\in 2^{\star}\times 2$, as
$$
\Sigma_{\mathbf{s},s}=
\left\{
\begin{array}{clll}
\{\alpha\},&\mbox{if }(\mathbf{s}, s)=(\lambda,0);\\
\{\beta\},&\mbox{if }(\mathbf{s},s)=((0,0,0),0);\\
\{\gamma\},&\mbox{if }(\mathbf{s},s)=((1,0),0);\\
\{\delta\},&\mbox{if }(\mathbf{s},s)=((1,1),0);\\
\varnothing&\mbox{otherwise;}
\end{array}\right.
$$
by $X$ be the $2$-sorted set defined as
$$
X_{s}=
\left\{
\begin{array}{lll}
\{a,b\},&\mbox{if }s=0;\\
\{x,y,z\},&\mbox{if }s=1;\\
\end{array}\right.
$$
and by $\mathcal{A}$ be the $2$-sorted subset of $\mathrm{T}_{\Sigma}(X)^{2}$ defined as
\[
\mathcal{A}_{s}=
\left\{
\begin{array}{ll}
\left\lbrace
\begin{array}{l}
\mathfrak{p}_{0}=(\alpha, \delta(x,z)),\\
\mathfrak{p}_{1}=(b,a),\\
\mathfrak{p}_{2}=(b,b),\\
\mathfrak{p}_{3}=(b,\delta(x,y)),\\
\mathfrak{p}_{4}=(\gamma(z,a),a),\\
\mathfrak{p}_{5}=(\delta(x,z), \beta(\gamma(z, a), b, a)),\\
\mathfrak{p}_{6}=(\beta(a,a,a),\delta(x,y)),\\
\mathfrak{p}_{7}=(\beta(a, \delta(z,x), a), \alpha)
\end{array}
\right\rbrace,
&\mbox{if } s=0;\\
\\
\left\lbrace
\begin{array}{l}
\mathfrak{q}_{0}=(x,z),\\
\mathfrak{q}_{1}=(y,x)
\end{array}
\right\rbrace,
&\mbox{if } s=1.
\end{array}\right.
\]

For the $2$-sorted specification $\boldsymbol{\mathcal{A}}=((\Sigma, X), \mathcal{A})$, let $\mathfrak{P}$ be the path in $\mathrm{Pth}_{\boldsymbol{\mathcal{A}},0}$ defined as the following sequence of steps
\allowdisplaybreaks
\begin{align*}
\mathfrak{P}\colon
\beta(\gamma(z, a),b,b)&
\xymatrix@C=100pt{\ar[r]^-{\text{\Small{$\left(\mathfrak{p}_{1}, 
\beta(\gamma(z, a),b,\underline{\quad})
\right)$}}}&}
\beta(\gamma(z, a), b, a)
\\&
\xymatrix@C=100pt{\ar[r]^-{\text{\Small{$\left(\mathfrak{p}_{2}, 
\beta(\gamma(z, a),\underline{\quad},a)
\right)$}}}&}
\beta(\gamma(z, a), b, a)
\\&
\xymatrix@C=100pt{\ar[r]^-{\text{\Small{$\left(\mathfrak{p}_{4}, 
\beta(\underline{\quad},b,a)
\right)$}}}&}
\beta(a, b, a)
\\&
\xymatrix@C=100pt{\ar[r]^-{\text{\Small{$\left(\mathfrak{p}_{1}, 
\beta(a,\underline{\quad},a)
\right)$}}}&}
\beta(a,a,a)
\\&
\xymatrix@C=100pt{\ar[r]^-{\text{\Small{$(\mathfrak{p}_{6}, 
\underline{\quad}
)$}}}&}
\delta(x,y)
\\&
\xymatrix@C=100pt{\ar[r]^-{\text{\Small{$(\mathfrak{q}_{1}, 
\delta(x,\underline{\quad})
)$}}}&}
\delta(x,x)
\\&
\xymatrix@C=100pt{\ar[r]^-{\text{\Small{$(\mathfrak{q}_{0}, 
\delta(x,\underline{\quad})
)$}}}&}
\delta(x,z)
\\&
\xymatrix@C=100pt{\ar[r]^-{\text{\Small{$(\mathfrak{p}_{5}, 
\underline{\quad}
)$}}}&}
\beta(\gamma(z, a), b, a)
\\&
\xymatrix@C=100pt{\ar[r]^-{\text{\Small{$(\mathfrak{p}_{3}, 
\beta(\gamma(z, a), \underline{\quad}, a)
)$}}}&}
\beta(\gamma(z, a), \delta(x,y), a)
\\&
\xymatrix@C=100pt{\ar[r]^-{\text{\Small{$(\mathfrak{q}_{1}, 
\beta(\gamma(z, a), \delta(x,\underline{\quad}), a)
)$}}}&}
\beta(\gamma(z, a), \delta(x,x), a)
\\&
\xymatrix@C=100pt{\ar[r]^-{\text{\Small{$(\mathfrak{q}_{0}, 
\beta(\gamma(z, a), \delta(\underline{\quad},x), a)
)$}}}&}
\beta(\gamma(z, a), \delta(z,x), a)
\\&
\xymatrix@C=100pt{\ar[r]^-{\text{\Small{$(\mathfrak{p}_{4}, 
\beta(\underline{\quad}, \delta(z,x), a)
)$}}}&}
\beta(a, \delta(z,x), a)
\\&
\xymatrix@C=100pt{\ar[r]^-{\text{\Small{$(\mathfrak{p}_{7}, 
\underline{\quad}
)$}}}&}
\alpha
\\&
\xymatrix@C=100pt{\ar[r]^-{\text{\Small{$(\mathfrak{p}_{0}, 
\underline{\quad}
)$}}}&}
\delta(x,z).
\end{align*}

By definition, $\mathfrak{P}$ is a $(00000110011000)$-path in $\boldsymbol{\mathcal{A}}$ of sort $0$ from $\beta(\gamma(z, a),b,b)$ to $\delta(x,z)$ and length $14$.

The following subpaths of $\mathfrak{P}$ are $0$-loops.
\[
\begin{array}{llll}
\mathfrak{P}^{1,1}&\colon
\beta(\gamma(z, a), b, a)&
\xymatrix@C=28pt{\ar[r]^-{}&}&
\beta(\gamma(z, a), b, a),
\\
\mathfrak{P}^{1,7}&\colon
\beta(\gamma(z, a), b, a)&
\xymatrix@C=28pt{\ar[r]^-{}&}&
\beta(\gamma(z, a), b, a),
\\
\mathfrak{P}^{2,7}&\colon
\beta(\gamma(z, a), b, a)&
\xymatrix@C=28pt{\ar[r]^-{}&}&
\beta(\gamma(z, a), b, a),
\\
\mathfrak{P}^{7,13}&\colon
\delta(x,z)&
\xymatrix@C=28pt{\ar[r]^-{}&}&
\delta(x,z).
\end{array}
\]

The following one-step subpaths of $\mathfrak{P}$ are all echelons.
\[
\begin{array}{llll}
\mathfrak{P}^{4,4}&\colon
\beta(a,a,a)&
\xymatrix@C=28pt{\ar[r]^-{}&}&
\delta(x,y),
\\
\mathfrak{P}^{7,7}&\colon
\delta(x,z)&
\xymatrix@C=28pt{\ar[r]^-{}&}&
\beta(\gamma(z, x), b,a),
\\
\mathfrak{P}^{12,12}&\colon
\beta(a,\delta(z,x),a)&
\xymatrix@C=28pt{\ar[r]^-{}&}&
\alpha,
\\
\mathfrak{P}^{13,13}&\colon
\alpha&
\xymatrix@C=28pt{\ar[r]^-{}&}&
\delta(x,z).
\end{array}
\]

Moreover, the following subpaths of $\mathfrak{P}$ are echelonless paths. 
\[
\begin{array}{llll}
\mathfrak{P}^{0,3}&\colon
\beta(\gamma(z, a), b,b)&
\xymatrix@C=28pt{\ar[r]^-{}&}&
\beta(a,a,a),
\\
\mathfrak{P}^{5,6}&\colon
\delta(x,y)&
\xymatrix@C=28pt{\ar[r]^-{}&}&
\delta(x,z),
\\
\mathfrak{P}^{8,11}&\colon
\beta(\gamma(z, a),b,a)&
\xymatrix@C=28pt{\ar[r]^-{}&}&
\beta(a,\delta(x,z), a).
\end{array}
\]
Let us note that, by Lemma~\ref{LPthHeadCt}, the just mentioned subpaths are head-constant paths.

For the echelonless subpath $\mathfrak{P}^{0,3}$ of $\mathfrak{P}$, the path extraction algorithm, from Lemma~\ref{LPthExtract}, applied to it retrieves the following family of paths.
\begin{center}
\begin{tikzpicture}[
fletxa/.style={thick, ->}]
\tikzstyle{every state}=[very thick, draw=blue!50,fill=blue!20, inner sep=2pt,minimum size=20pt]
\node[] (0P) at (0,0) [] {$P_{0}$};
\node[] (0=) at (.7,0) [] {$=$};
\node[] (0s) at (1.2,0) [] {$\beta($};
\node[] (00) at (2.2,0) [] {$\color{white}{P^{i}_{\bb{\mathbf{s}}-1}}$};
\node[] () at (2.2,0) [] {$\gamma(z, a)$};
\node[] (0c0) at (3.4,0) [] {$,$};
\node[] (0j) at (4.5,0) [] {$\color{white}{P^{i}_{\bb{\mathbf{s}}-1}}$};
\node[] () at (4.5,0) [] {$b$};
\node[] (0c1) at (5.6,0) [] {$,$};
\node[] (0f) at (6.7,0) [] {$\color{white}{P^{i}_{\bb{\mathbf{s}}-1}}$};
\node[] () at (6.7,0) [] {$b$};
\node[] (0c0) at (7.6,0) [] {$)$};

\node[] (E0label) at (-1.3,-.75) [] {$\scriptstyle (\mathfrak{p}_{1}, \beta(\gamma(z,a),b,\underline{\quad}))$};
\node[] (E0) at (0,-.75) [] {\rotatebox{-90}{$\mor$}};
\node[] (E00) at (2.2,-.75) [] {\rotatebox{-90}{$\,=$}};
\node[] (E0j) at (4.5,-.75) [] {\rotatebox{-90}{$\,=$}};
\node[] (E0f) at (6.7,-.75) [] {\rotatebox{-90}{$\mor$}};
\node[] (E0flabel) at (7.5,-.75) [] {$\scriptstyle (\mathfrak{p}_{1},\underline{\quad})$};

\node[] (1P) at (0,-1.5) [] {$P_{1}$};
\node[] (1=) at (.7,-1.5) [] {$=$};
\node[] (1s) at (1.2,-1.5) [] {$\beta($};
\node[] (10) at (2.2,-1.5) [] {$\color{white}P^{i}_{\bb{\mathbf{s}}-1}$};
\node[] () at (2.2,-1.5) [] {$\gamma(z, a)$};
\node[] (1c0) at (3.4,-1.5) [] {$, $};
\node[] (1j) at (4.5,-1.5) [] {$\color{white}P^{i}_{\bb{\mathbf{s}}-1}$};
\node[] () at (4.5,-1.5) [] {$b$};
\node[] (1c1) at (5.6,-1.5) [] {$, $};
\node[] (1f) at (6.7,-1.5) [] {$\color{white}P^{i}_{\bb{\mathbf{s}}-1}$};
\node[] () at (6.7,-1.5) [] {$a$};
\node[] (1c0) at (7.6,-1.5) [] {$)$};

\node[] (E1label) at (-1.3,-2.25) [] {$\scriptstyle (\mathfrak{p}_{2}, \beta(\gamma(z,a),\underline{\quad},a))$};
\node[] (E1) at (0,-2.25) [] {\rotatebox{-90}{$\mor$}};
\node[] (E10) at (2.2,-2.25) [] {\rotatebox{-90}{$\,=$}};
\node[] (E1j) at (4.5,-2.25) [] {\rotatebox{-90}{$\mor$}};
\node[] (E1jlabel) at (5.3,-2.25) [] {$\scriptstyle (\mathfrak{p}_{2}, \underline{\quad})$};
\node[] (E1f) at (6.7,-2.25) [] {\rotatebox{-90}{$\,=$}};

\node[] (2P) at (0,-3) [] {$P_{2}$};
\node[] (2=) at (.7,-3) [] {$=$};
\node[] (2s) at (1.2,-3) [] {$\beta($};
\node[] (20) at (2.2,-3) [] {$\color{white}P^{i}_{\bb{\mathbf{s}}-1}$};
\node[] () at (2.2,-3) [] {$\gamma(z, a)$};
\node[] (2c0) at (3.4,-3) [] {$, $};
\node[] (2j) at (4.5,-3) [] {$\color{white}P^{i}_{\bb{\mathbf{s}}-1}$};
\node[] () at (4.5,-3) [] {$b$};
\node[] (2c1) at (5.6,-3) [] {$, $};
\node[] (2f) at (6.7,-3) [] {$\color{white}P^{i}_{\bb{\mathbf{s}}-1}$};
\node[] () at (6.7,-3) [] {$a$};
\node[] (2c0) at (7.6,-3) [] {$)$};

\node[] (E2label) at (-1,-3.75) [] {$\scriptstyle (\mathfrak{p}_{4}, \beta(\underline{\quad},b,a))$};
\node[] (E2) at (0,-3.75) [] {\rotatebox{-90}{$\mor$}};
\node[] (E20) at (2.2,-3.75) [] {\rotatebox{-90}{$\mor$}};
\node[] (E20label) at (2.85,-3.75) [] {$\scriptstyle (\mathfrak{p}_{4}, \underline{\quad})$};
\node[] (E2j) at (4.5,-3.75) [] {\rotatebox{-90}{$\,=$}};
\node[] (E2f) at (6.7,-3.75) [] {\rotatebox{-90}{$\,=$}};

\node[] (3P) at (0,-4.5) [] {$P_{3}$};
\node[] (3=) at (.7,-4.5) [] {$=$};
\node[] (3s) at (1.2,-4.5) [] {$\beta($};
\node[] (30) at (2.2,-4.5) [] {$\color{white}P^{i}_{\bb{\mathbf{s}}-1}$};
\node[] () at (2.2,-4.5) [] {$ a$};
\node[] (3c0) at (3.4,-4.5) [] {$, $};
\node[] (3j) at (4.5,-4.5) [] {$\color{white}P^{i}_{\bb{\mathbf{s}}-1}$};
\node[] () at (4.5,-4.5) [] {$b$};
\node[] (3c1) at (5.6,-4.5) [] {$, $};
\node[] (3f) at (6.7,-4.5) [] {$\color{white}P^{i}_{\bb{\mathbf{s}}-1}$};
\node[] () at (6.7,-4.5) [] {$a$};
\node[] (3c0) at (7.6,-4.5) [] {$)$};

\node[] (E3label) at (-1,-5.25) [] {$\scriptstyle (\mathfrak{p}_{1},\beta(a,\underline{\quad},a))$};
\node[] (E3) at (0,-5.25) [] {\rotatebox{-90}{$\mor$}};
\node[] (E30) at (2.2,-5.25) [] {\rotatebox{-90}{$\,=$}};
\node[] (E3jlabel) at (5.15,-5.25) [] {$\scriptstyle (\mathfrak{p}_{1},\underline{\quad})$};
\node[] (E3j) at (4.5,-5.25) [] {\rotatebox{-90}{$\mor$}};
\node[] (E3f) at (6.7,-5.25) [] {\rotatebox{-90}{$\,=$}};

\node[] (4P) at (0,-6) [] {$P_{4}$};
\node[] (4=) at (.7,-6) [] {$=$};
\node[] (4s) at (1.2,-6) [] {$\beta($};
\node[] (40) at (2.2,-6) [] {$\color{white}P^{i}_{\bb{\mathbf{s}}-1}$};
\node[] () at (2.2,-6) [] {$ a$};
\node[] (4c0) at (3.4,-6) [] {$, $};
\node[] (4j) at (4.5,-6) [] {$\color{white}P^{i}_{\bb{\mathbf{s}}-1}$};
\node[] () at (4.5,-6) [] {$a$};
\node[] (4c1) at (5.6,-6) [] {$, $};
\node[] (4f) at (6.7,-6) [] {$\color{white}P^{i}_{\bb{\mathbf{s}}-1}$};
\node[] () at (6.7,-6) [] {$a$};
\node[] (4c0) at (7.6,-6) [] {$)$};

\tikzset{encercla/.style={draw=black, line width=.5pt, inner sep=0pt, rectangle, rounded corners}};
\node [encercla,  fit=(00) (10) (20) ] {} ;
\node [encercla,  fit=(30) (40) ] {} ;
\node [encercla,  fit=(0j) (1j) ] {} ;
\node [encercla,  fit=(2j) (3j) ] {} ;
\node [encercla,  fit=(4j) ] {} ;
\node [encercla,  fit=(0f) ] {} ;
\node [encercla,  fit=(1f) (2f) (3f) (4f) ] {} ;
\end{tikzpicture}
\end{center}
The paths that we can extract from $\mathfrak{P}^{0,3}$ are, thus, the components of the family of paths 
$(\mathfrak{P}_{j})_{j\in 3}$ in $\mathrm{Pth}_{\boldsymbol{\mathcal{A}},(0,0,0)}$ where 
\[
\begin{array}{lllll}
\mathfrak{P}_{0}\colon \gamma(z, a)&
\xymatrix@C=40pt{\ar[r]^-{\text{\Small{$(\mathfrak{p}_{4},\underline{\quad})$}}}&}&
a;
\\
\mathfrak{P}_{1}\colon b&
\xymatrix@C=40pt{\ar[r]^-{\text{\Small{$(\mathfrak{p}_{2},\underline{\quad})$}}}&}&
b&
\xymatrix@C=40pt{\ar[r]^-{\text{\Small{$(\mathfrak{p}_{1},\underline{\quad})$}}}&}&
a;
\\
\mathfrak{P}_{2}\colon  b&
\xymatrix@C=40pt{\ar[r]^-{\text{\Small{$(\mathfrak{p}_{1},\underline{\quad})$}}}&}&
a.
\end{array}
\]
Let us note that (1) $\mathfrak{P}^{0,3}$ is a $(0000)$-path in $\boldsymbol{\mathcal{A}}$ of sort $0$, (2) $\mathfrak{P}_{0}$ is a $0$-path in $\boldsymbol{\mathcal{A}}$ of sort $0$, (3) $\mathfrak{P}_{1}$ is a $(00)$-path in $\boldsymbol{\mathcal{A}}$ of sort $0$, and (4) $\mathfrak{P}_{2}$ is a $0$-path  in $\boldsymbol{\mathcal{A}}$ of sort $0$. In particular, $\mathfrak{P}_{0}$ and $\mathfrak{P}_{2}$ are echelons, and $\mathfrak{P}_{1}$ is the composition of two echelons.

Consider the $(0000)$-path in $\boldsymbol{\mathcal{A}}$ of sort $0$ given by $\beta^{\mathbf{Pth}_{\boldsymbol{\mathcal{A}}}}((\mathfrak{P}_{j})_{j\in 3})$, i.e., the result of applying $\beta$ to the family of paths $(\mathfrak{P}_{j})_{j\in 3}$ in the $2$-sorted $\Sigma$-algebra $\mathbf{Pth}^{(0,1)}_{\boldsymbol{\mathcal{A}}}$. This path, by Proposition~\ref{PPthAlg}, is given by the following sequence
\allowdisplaybreaks
\begin{align*}
\beta^{\mathbf{Pth}_{\boldsymbol{\mathcal{A}}}}((\mathfrak{P}_{j})_{j\in 3})\colon
\beta(\gamma(z, a),b,b)&
\xymatrix@C=100pt{\ar[r]^-{\text{\Small{$\left(\mathfrak{p}_{4}, 
\beta(\underline{\quad},b,a)
\right)$}}}&}
\beta(a,b,b)
\\&
\xymatrix@C=100pt{\ar[r]^-{\text{\Small{$\left(\mathfrak{p}_{2}, 
\beta(a,\underline{\quad},a)
\right)$}}}&}
\beta(a,b,b)
\\&
\xymatrix@C=100pt{\ar[r]^-{\text{\Small{$\left(\mathfrak{p}_{1}, 
\beta(a,\underline{\quad},a)
\right)$}}}&}
\beta(a,a,b)
\\&
\xymatrix@C=100pt{\ar[r]^-{\text{\Small{$\left(\mathfrak{p}_{1}, 
\beta(a,a,\underline{\quad})
\right)$}}}&}
\beta(a,a,a).
\end{align*}

The paths that we can extract from $\beta^{\mathbf{Pth}_{\boldsymbol{\mathcal{A}}}}((\mathfrak{P}_{j})_{j\in 3})$ are, by Proposition~\ref{PRecov}, again, the components of the family of paths $(\mathfrak{P}_{j})_{j\in 3}$. However, the appearance is no longer chaotic. We see how the path extraction algorithm, from Lemma~\ref{LPthExtract}, applies to $\beta^{\mathbf{Pth}_{\boldsymbol{\mathcal{A}}}}((\mathfrak{P}_{j})_{j\in 3})$.
\begin{center}
\begin{tikzpicture}[
fletxa/.style={thick, ->}]
\tikzstyle{every state}=[very thick, draw=blue!50,fill=blue!20, inner sep=2pt,minimum size=20pt]
\node[] (0P) at (0,0) [] {$P_{0}$};
\node[] (0=) at (.7,0) [] {$=$};
\node[] (0s) at (1.2,0) [] {$\beta($};
\node[] (00) at (2.2,0) [] {$\color{white}{P^{i}_{\bb{\mathbf{s}}-1}}$};
\node[] () at (2.2,0) [] {$\gamma(z, a)$};
\node[] (0c0) at (3.4,0) [] {$,$};
\node[] (0j) at (4.5,0) [] {$\color{white}{P^{i}_{\bb{\mathbf{s}}-1}}$};
\node[] () at (4.5,0) [] {$b$};
\node[] (0c1) at (5.6,0) [] {$,$};
\node[] (0f) at (6.7,0) [] {$\color{white}{P^{i}_{\bb{\mathbf{s}}-1}}$};
\node[] () at (6.7,0) [] {$b$};
\node[] (0c0) at (7.6,0) [] {$)$};

\node[] (E0label) at (-1,-.75) [] {$\scriptstyle (\mathfrak{p}_{4}, \beta(\underline{\quad},b,b))$};
\node[] (E0) at (0,-.75) [] {\rotatebox{-90}{$\mor$}};
\node[] (E00) at (2.2,-.75) [] {\rotatebox{-90}{$\mor$}};
\node[] (E00label) at (2.85,-.75) [] {$\scriptstyle (\mathfrak{p}_{4},\underline{\quad})$};
\node[] (E0j) at (4.5,-.75) [] {\rotatebox{-90}{$\,=$}};
\node[] (E0f) at (6.7,-.75) [] {\rotatebox{-90}{$\,=$}};

\node[] (1P) at (0,-1.5) [] {$Q_{1}$};
\node[] (1=) at (.7,-1.5) [] {$=$};
\node[] (1s) at (1.2,-1.5) [] {$\beta($};
\node[] (10) at (2.2,-1.5) [] {$\color{white}P^{i}_{\bb{\mathbf{s}}-1}$};
\node[] () at (2.2,-1.5) [] {$a$};
\node[] (1c0) at (3.4,-1.5) [] {$, $};
\node[] (1j) at (4.5,-1.5) [] {$\color{white}P^{i}_{\bb{\mathbf{s}}-1}$};
\node[] () at (4.5,-1.5) [] {$b$};
\node[] (1c1) at (5.6,-1.5) [] {$, $};
\node[] (1f) at (6.7,-1.5) [] {$\color{white}P^{i}_{\bb{\mathbf{s}}-1}$};
\node[] () at (6.7,-1.5) [] {$b$};
\node[] (1c0) at (7.6,-1.5) [] {$)$};

\node[] (E1label) at (-1,-2.25) [] {$\scriptstyle (\mathfrak{p}_{2},\beta(a,\underline{\quad},a)$};
\node[] (E1) at (0,-2.25) [] {\rotatebox{-90}{$\mor$}};
\node[] (E10) at (2.2,-2.25) [] {\rotatebox{-90}{$\,=$}};
\node[] (E1j) at (4.5,-2.25) [] {\rotatebox{-90}{$\mor$}};
\node[] (E1jlabel) at (5.15,-2.25) [] {$\scriptstyle (\mathfrak{p}_{2},\underline{\quad})$};
\node[] (E1f) at (6.7,-2.25) [] {\rotatebox{-90}{$\,=$}};

\node[] (2P) at (0,-3) [] {$Q_{2}$};
\node[] (2=) at (.7,-3) [] {$=$};
\node[] (2s) at (1.2,-3) [] {$\beta($};
\node[] (20) at (2.2,-3) [] {$\color{white}P^{i}_{\bb{\mathbf{s}}-1}$};
\node[] () at (2.2,-3) [] {$a$};
\node[] (2c0) at (3.4,-3) [] {$, $};
\node[] (2j) at (4.5,-3) [] {$\color{white}P^{i}_{\bb{\mathbf{s}}-1}$};
\node[] () at (4.5,-3) [] {$b$};
\node[] (2c1) at (5.6,-3) [] {$, $};
\node[] (2f) at (6.7,-3) [] {$\color{white}P^{i}_{\bb{\mathbf{s}}-1}$};
\node[] () at (6.7,-3) [] {$b$};
\node[] (2c0) at (7.6,-3) [] {$)$};

\node[] (E2label) at (-1,-3.75) [] {$\scriptstyle (\mathfrak{p}_{1},\beta(a,\underline{\quad},a))$};
\node[] (E2) at (0,-3.75) [] {\rotatebox{-90}{$\mor$}};
\node[] (E20) at (2.2,-3.75) [] {\rotatebox{-90}{$\,=$}};
\node[] (E2j) at (4.5,-3.75) [] {\rotatebox{-90}{$\mor$}};
\node[] (E2jlabel) at (5.15,-3.75) [] {$\scriptstyle (\mathfrak{p}_{1},\underline{\quad})$};
\node[] (E2f) at (6.7,-3.75) [] {\rotatebox{-90}{$\,=$}};

\node[] (3P) at (0,-4.5) [] {$Q_{3}$};
\node[] (3=) at (.7,-4.5) [] {$=$};
\node[] (3s) at (1.2,-4.5) [] {$\beta($};
\node[] (30) at (2.2,-4.5) [] {$\color{white}P^{i}_{\bb{\mathbf{s}}-1}$};
\node[] () at (2.2,-4.5) [] {$ a$};
\node[] (3c0) at (3.4,-4.5) [] {$, $};
\node[] (3j) at (4.5,-4.5) [] {$\color{white}P^{i}_{\bb{\mathbf{s}}-1}$};
\node[] () at (4.5,-4.5) [] {$a$};
\node[] (3c1) at (5.6,-4.5) [] {$, $};
\node[] (3f) at (6.7,-4.5) [] {$\color{white}P^{i}_{\bb{\mathbf{s}}-1}$};
\node[] () at (6.7,-4.5) [] {$b$};
\node[] (3c0) at (7.6,-4.5) [] {$)$};

\node[] (E3label) at (-1,-5.25) [] {$\scriptstyle (\mathfrak{p}_{1},\beta(a,a,\underline{\quad}))$};
\node[] (E3) at (0,-5.25) [] {\rotatebox{-90}{$\mor$}};
\node[] (E30) at (2.2,-5.25) [] {\rotatebox{-90}{$\,=$}};
\node[] (E3j) at (4.5,-5.25) [] {\rotatebox{-90}{$\,=$}};
\node[] (E3flabel) at (7.35,-5.25) [] {$\scriptstyle (\mathfrak{p}_{1},\underline{\quad})$};
\node[] (E3f) at (6.7,-5.25) [] {\rotatebox{-90}{$\mor$}};

\node[] (4P) at (0,-6) [] {$P_{4}$};
\node[] (4=) at (.7,-6) [] {$=$};
\node[] (4s) at (1.2,-6) [] {$\beta($};
\node[] (40) at (2.2,-6) [] {$\color{white}P^{i}_{\bb{\mathbf{s}}-1}$};
\node[] () at (2.2,-6) [] {$ a$};
\node[] (4c0) at (3.4,-6) [] {$, $};
\node[] (4j) at (4.5,-6) [] {$\color{white}P^{i}_{\bb{\mathbf{s}}-1}$};
\node[] () at (4.5,-6) [] {$a$};
\node[] (4c1) at (5.6,-6) [] {$, $};
\node[] (4f) at (6.7,-6) [] {$\color{white}P^{i}_{\bb{\mathbf{s}}-1}$};
\node[] () at (6.7,-6) [] {$a$};
\node[] (4c0) at (7.6,-6) [] {$)$};

\tikzset{encercla/.style={draw=black, line width=.5pt, inner sep=0pt, rectangle, rounded corners}};
\node [encercla,  fit=(00)  ] {} ;
\node [encercla,  fit=(10) (20) (30) (40)  ] {} ;
\node [encercla,  fit=(0j) (1j)  ] {} ;
\node [encercla,  fit=(2j)  ] {} ;
\node [encercla,  fit=(3j) (4j)  ] {} ;
\node [encercla,  fit=(0f) (1f) (2f) (3f)] {} ;
\node [encercla,  fit=(4f)  ] {} ;
\end{tikzpicture}
\end{center}
The path $\beta^{\mathbf{Pth}_{\boldsymbol{\mathcal{A}}}}((\mathfrak{P}_{j})_{j\in 3})$ is a path where all rewrite rules are performed following a leftmost innermost derivation strategy. Although  the intermediate terms may have been changed, this path has the same $(0,1)$-source and $(0,1)$-target as $\mathfrak{P}^{0,3}$.

For the echelonless subpath $\mathfrak{P}^{5,6}$ from $\mathfrak{P}$, the path extraction algorithm, from Lemma~\ref{LPthExtract}, applied to it retrieves the following family of paths.
\begin{center}
\begin{tikzpicture}[
fletxa/.style={thick, ->}]
\tikzstyle{every state}=[very thick, draw=blue!50,fill=blue!20, inner sep=2pt,minimum size=20pt]
\node[] (0P) at (0,0) [] {$P_{5}$};
\node[] (0=) at (.7,0) [] {$=$};
\node[] (0s) at (1.2,0) [] {$\delta($};
\node[] (00) at (2.2,0) [] {$\color{white}{P^{i}_{\bb{\mathbf{s}}-1}}$};
\node[] () at (2.2,0) [] {$x$};
\node[] (0c0) at (3.4,0) [] {$,$};
\node[] (0j) at (4.5,0) [] {$\color{white}{P^{i}_{\bb{\mathbf{s}}-1}}$};
\node[] () at (4.5,0) [] {$y$};
\node[] () at (5.4,0) [] {$)$};

\node[] (E0label) at (-1,-.75) [] {$\scriptstyle (\mathfrak{q}_{1},\delta(x,\underline{\quad}))$};
\node[] (E0) at (-0,-.75) [] {\rotatebox{-90}{$\mor$}};
\node[] (E00) at (2.2,-.75) [] {\rotatebox{-90}{$\,=$}};
\node[] (E0f) at (4.5,-.75) [] {\rotatebox{-90}{$\mor$}};
\node[] (E0flabel) at (5.15,-.75) [] {$\scriptstyle (\mathfrak{q}_{1},\underline{\quad})$};

\node[] (1P) at (0,-1.5) [] {$P_{6}$};
\node[] (1=) at (.7,-1.5) [] {$=$};
\node[] (1s) at (1.2,-1.5) [] {$\delta($};
\node[] (10) at (2.2,-1.5) [] {$\color{white}{P^{i}_{\bb{\mathbf{s}}-1}}$};
\node[] () at (2.2,-1.5) [] {$x$};
\node[] (1c0) at (3.4,-1.5) [] {$,$};
\node[] (1j) at (4.5,-1.5) [] {$\color{white}{P^{i}_{\bb{\mathbf{s}}-1}}$};
\node[] () at (4.5,-1.5) [] {$x$};
\node[] () at (5.4,-1.5) [] {$)$};

\node[] (E1label) at (-1,-2.25) [] {$\scriptstyle (\mathfrak{q}_{0},\delta(x,\underline{\quad}))$};
\node[] (E1) at (-0,-2.25) [] {\rotatebox{-90}{$\mor$}};
\node[] (E10) at (2.2,-2.25) [] {\rotatebox{-90}{$\,=$}};
\node[] (E1f) at (4.5,-2.25) [] {\rotatebox{-90}{$\mor$}};
\node[] (E1flabel) at (5.15,-2.25) [] {$\scriptstyle (\mathfrak{q}_{0},\underline{\quad})$};

\node[] (2P) at (0,-3) [] {$P_{7}$};
\node[] (2=) at (.7,-3) [] {$=$};
\node[] (2s) at (1.2,-3) [] {$\delta($};
\node[] (20) at (2.2,-3) [] {$\color{white}{P^{i}_{\bb{\mathbf{s}}-1}}$};
\node[] () at (2.2,-3) [] {$x$};
\node[] (2c0) at (3.4,-3) [] {$,$};
\node[] (2j) at (4.5,-3) [] {$\color{white}{P^{i}_{\bb{\mathbf{s}}-1}}$};
\node[] () at (4.5,-3) [] {$z$};
\node[] () at (5.4,-3) [] {$)$};

\tikzset{encercla/.style={draw=black, line width=.5pt, inner sep=0pt, rectangle, rounded corners}};
\node [encercla,  fit=(00) (10) (20) ] {} ;
\node [encercla,  fit=(0j) ] {} ;
\node [encercla,  fit=(1j) ] {} ;
\node [encercla,  fit=(2j) ] {} ;
\end{tikzpicture}
\end{center}
The paths that we can extract from $\mathfrak{P}^{5,6}$ are, thus, the components of the family of paths $(\mathfrak{Q}_{j})_{j\in 2}$ in $\mathrm{Pth}_{\boldsymbol{\mathcal{A}},(1,1)}$ where
\[
\begin{array}{lllll}
\mathfrak{Q}_{0}\colon x
\\
\mathfrak{Q}_{1}\colon y&
\xymatrix@C=50pt{\ar[r]^-{\text{\Small{$(\mathfrak{q}_{1},\underline{\quad})$}}}&}&
x&
\xymatrix@C=50pt{\ar[r]^-{\text{\Small{$(\mathfrak{q}_{0},\underline{\quad})$}}}&}&
z.
\end{array}
\]

Let us note that (1) $\mathfrak{P}^{5,6}$ is a $(11)$-path in $\boldsymbol{\mathcal{A}}$ of sort $0$, (2) $\mathfrak{Q}_{0}$ is a $\lambda$-path in $\boldsymbol{\mathcal{A}}$ of sort $1$, and (3) $\mathfrak{Q}_{1}$ is a $(11)$-path in $\boldsymbol{\mathcal{A}}$ of sort $1$. In particular, $\mathfrak{Q}_{0}$ is the $(1,0)$-identity path on $x$ and $\mathfrak{Q}_{1}$ is the composition of two echelons.

Consider the $(11)$-path in $\boldsymbol{\mathcal{A}}$ of sort $0$ given by $\delta^{\mathbf{Pth}_{\boldsymbol{\mathcal{A}}}}(\mathfrak{Q}_{0}, \mathfrak{Q}_{1})$, i.e., the result of applying $\delta$ to the family of paths $(\mathfrak{Q}_{j})_{j\in 2}$ in the $2$-sorted $\Sigma$-algebra $\mathbf{Pth}^{(0,1)}_{\boldsymbol{\mathcal{A}}}$. This path, by Proposition~\ref{PPthAlg}, is given by the following sequence. 
\[
\begin{array}{lllll}
\delta^{\mathbf{Pth}_{\boldsymbol{\mathcal{A}}}}(\mathfrak{Q}_{0}, \mathfrak{Q}_{1})\colon \delta(x,y)&
\xymatrix@C=60pt{\ar[r]^-{\text{\Small{$(\mathfrak{q}_{1},\delta(x,\underline{\quad}))$}}}&}&
\delta(x,x)&
\xymatrix@C=60pt{\ar[r]^-{\text{\Small{$(\mathfrak{q}_{0},\delta(x,\underline{\quad}))$}}}&}&
\delta(x,z).
\end{array}
\]
It follows that $\mathfrak{P}^{5,6}=\delta^{\mathbf{Pth}_{\boldsymbol{\mathcal{A}}}}(\mathfrak{Q}_{0}, \mathfrak{Q}_{1})$, i.e., that extracting paths from $\mathfrak{P}^{5,6}$ and, afterwards, considering the result of the $\delta$ operation in $\mathbf{Pth}^{(0,1)}_{\boldsymbol{\mathcal{A}}}$ on the family $(\mathfrak{Q}_{0}, \mathfrak{Q}_{1})$ does not alter the original path. This is so, because, in the original path, $\mathfrak{P}^{5,6}$, the derivation process followed a leftmost innermost derivation strategy. Let us note that the positions marking where to perform the transformations in a path, following a leftmost innermost transformation strategy, are not necessarily equal to $0$.

For the echelonless subpath $\mathfrak{P}^{8,11}$ from $\mathfrak{P}$, the path extraction algorithm, from Lemma~\ref{LPthExtract}, applied to it retrieves the following family of paths.
\begin{center}
\begin{tikzpicture}[
fletxa/.style={thick, ->}]
\tikzstyle{every state}=[very thick, draw=blue!50,fill=blue!20, inner sep=2pt,minimum size=20pt]
\node[] (0P) at (0,0) [] {$P_{8}$};
\node[] (0=) at (.7,0) [] {$=$};
\node[] (0s) at (1.2,0) [] {$\beta($};
\node[] (00) at (2.2,0) [] {$\color{white}{P^{i}_{\bb{\mathbf{s}}-1}}$};
\node[] () at (2.2,0) [] {$\gamma(z, a)$};
\node[] (0c0) at (3.4,0) [] {$,$};
\node[] (0j) at (4.5,0) [] {$\color{white}{P^{i}_{\bb{\mathbf{s}}-1}}$};
\node[] () at (4.5,0) [] {$b$};
\node[] (0c1) at (5.6,0) [] {$,$};
\node[] (0f) at (6.7,0) [] {$\color{white}{P^{i}_{\bb{\mathbf{s}}-1}}$};
\node[] () at (6.7,0) [] {$a$};
\node[] (0c0) at (7.6,0) [] {$)$};

\node[] (E0label) at (-1.3,-.75) [] {$\scriptstyle (\mathfrak{p}_{3}, \beta(\gamma(z,a),\underline{\quad},a))$};
\node[] (E0) at (0,-.75) [] {\rotatebox{-90}{$\mor$}};
\node[] (E00) at (2.2,-.75) [] {\rotatebox{-90}{$\,=$}};
\node[] (E0j) at (4.5,-.75) [] {\rotatebox{-90}{$\mor$}};
\node[] (E0flabel) at (5.15,-.75) [] {$\scriptstyle (\mathfrak{p}_{3},\underline{\quad})$};
\node[] (E0f) at (6.7,-.75) [] {\rotatebox{-90}{$\,=$}};

\node[] (1P) at (0,-1.5) [] {$P_{9}$};
\node[] (1=) at (.7,-1.5) [] {$=$};
\node[] (1s) at (1.2,-1.5) [] {$\beta($};
\node[] (10) at (2.2,-1.5) [] {$\color{white}P^{i}_{\bb{\mathbf{s}}-1}$};
\node[] () at (2.2,-1.5) [] {$\gamma(z, a)$};
\node[] (1c0) at (3.4,-1.5) [] {$, $};
\node[] (1j) at (4.5,-1.5) [] {$\color{white}P^{i}_{\bb{\mathbf{s}}-1}$};
\node[] () at (4.5,-1.5) [] {$\delta(x,y)$};
\node[] (1c1) at (5.6,-1.5) [] {$, $};
\node[] (1f) at (6.7,-1.5) [] {$\color{white}P^{i}_{\bb{\mathbf{s}}-1}$};
\node[] () at (6.7,-1.5) [] {$a$};
\node[] (1c0) at (7.6,-1.5) [] {$)$};

\node[] (E1label) at (-1.65,-2.25) [] {$\scriptstyle (\mathfrak{q}_{1}, \beta(\gamma(z,a),\delta(x,\underline{\quad}),a))$};
\node[] (E1) at (0,-2.25) [] {\rotatebox{-90}{$\mor$}};
\node[] (E10) at (2.2,-2.25) [] {\rotatebox{-90}{$\,=$}};
\node[] (E1j) at (4.5,-2.25) [] {\rotatebox{-90}{$\mor$}};
\node[] (E1jlabel) at (5.35,-2.25) [] {$\scriptstyle (\mathfrak{q}_{1}, \delta(x,\underline{\;\;\;}))$};
\node[] (E1f) at (6.7,-2.25) [] {\rotatebox{-90}{$\,=$}};

\node[] (2P) at (0,-3) [] {$P_{10}$};
\node[] (2=) at (.7,-3) [] {$=$};
\node[] (2s) at (1.2,-3) [] {$\beta($};
\node[] (20) at (2.2,-3) [] {$\color{white}P^{i}_{\bb{\mathbf{s}}-1}$};
\node[] () at (2.2,-3) [] {$\gamma(z, a)$};
\node[] (2c0) at (3.4,-3) [] {$, $};
\node[] (2j) at (4.5,-3) [] {$\color{white}P^{i}_{\bb{\mathbf{s}}-1}$};
\node[] () at (4.5,-3) [] {$\delta(x,x)$};
\node[] (2c1) at (5.6,-3) [] {$, $};
\node[] (2f) at (6.7,-3) [] {$\color{white}P^{i}_{\bb{\mathbf{s}}-1}$};
\node[] () at (6.7,-3) [] {$a$};
\node[] (2c0) at (7.6,-3) [] {$)$};

\node[] (E2label) at (-1.65,-3.75) [] {$\scriptstyle (\mathfrak{q}_{0}, \beta(\gamma(z,a),\delta(\underline{\quad},x),a))$};
\node[] (E2) at (0,-3.75) [] {\rotatebox{-90}{$\mor$}};
\node[] (E20) at (2.2,-3.75) [] {\rotatebox{-90}{$\,=$}};
\node[] (E2j) at (4.5,-3.75) [] {\rotatebox{-90}{$\mor$}};
\node[] (E2jlabel) at (5.35,-3.75) [] {$\scriptstyle (\mathfrak{q}_{0}, \delta(\underline{\;\;\;},x))$};
\node[] (E2f) at (6.7,-3.75) [] {\rotatebox{-90}{$\,=$}};

\node[] (3P) at (0,-4.5) [] {$P_{11}$};
\node[] (3=) at (.7,-4.5) [] {$=$};
\node[] (3s) at (1.2,-4.5) [] {$\beta($};
\node[] (30) at (2.2,-4.5) [] {$\color{white}P^{i}_{\bb{\mathbf{s}}-1}$};
\node[] () at (2.2,-4.5) [] {$ \gamma(z, a)$};
\node[] (3c0) at (3.4,-4.5) [] {$, $};
\node[] (3j) at (4.5,-4.5) [] {$\color{white}P^{i}_{\bb{\mathbf{s}}-1}$};
\node[] () at (4.5,-4.5) [] {$\delta(z,x)$};
\node[] (3c1) at (5.6,-4.5) [] {$, $};
\node[] (3f) at (6.7,-4.5) [] {$\color{white}P^{i}_{\bb{\mathbf{s}}-1}$};
\node[] () at (6.7,-4.5) [] {$a$};
\node[] (3c0) at (7.6,-4.5) [] {$)$};

\node[] (E3label) at (-1.5,-5.25) [] {$\scriptstyle (\mathfrak{p}_{4},\beta(\underline{\quad},\delta(z,x),a))$};
\node[] (E3) at (0,-5.25) [] {\rotatebox{-90}{$\mor$}};
\node[] (E30label) at (2.85,-5.25) [] {$\scriptstyle (\mathfrak{p}_{4}, \underline{\quad})$};
\node[] (E30) at (2.2,-5.25) [] {\rotatebox{-90}{$\mor$}};
\node[] (E3j) at (4.5,-5.25) [] {\rotatebox{-90}{$\,=$}};
\node[] (E3f) at (6.7,-5.25) [] {\rotatebox{-90}{$\,=$}};

\node[] (4P) at (0,-6) [] {$P_{12}$};
\node[] (4=) at (.7,-6) [] {$=$};
\node[] (4s) at (1.2,-6) [] {$\beta($};
\node[] (40) at (2.2,-6) [] {$\color{white}P^{i}_{\bb{\mathbf{s}}-1}$};
\node[] () at (2.2,-6) [] {$ a$};
\node[] (4c0) at (3.4,-6) [] {$, $};
\node[] (4j) at (4.5,-6) [] {$\color{white}P^{i}_{\bb{\mathbf{s}}-1}$};
\node[] () at (4.5,-6) [] {$\delta(z,x)$};
\node[] (4c1) at (5.6,-6) [] {$, $};
\node[] (4f) at (6.7,-6) [] {$\color{white}P^{i}_{\bb{\mathbf{s}}-1}$};
\node[] () at (6.7,-6) [] {$a$};
\node[] (4c0) at (7.6,-6) [] {$)$};

\tikzset{encercla/.style={draw=black, line width=.5pt, inner sep=0pt, rectangle, rounded corners}};
\node [encercla,  fit=(00) (10) (20) (30) ] {} ;
\node [encercla,  fit=(40)  ] {} ;
\node [encercla,  fit=(0f) (1f) (2f) (3f) (4f) ] {} ;
\node [encercla,  fit=(0j)  ] {} ;
\node [encercla,  fit=(1j)  ] {} ;
\node [encercla,  fit=(2j)  ] {} ;
\node [encercla,  fit=(3j) (4j)  ] {} ;
\end{tikzpicture}
\end{center}
The paths that we can extract from $\mathfrak{P}^{8,11}$ are, thus, the components of the family of paths $(\mathfrak{R}_{j})_{j\in 3}$ in $\mathrm{Pth}_{\boldsymbol{\mathcal{A}},000}$ where
\[
\allowdisplaybreaks
\begin{array}{lllllll}
\mathfrak{R}_{0}\colon \gamma(z, a)&
\xymatrix@C=40pt{\ar[r]^-{\text{\Small{$(\mathfrak{p}_{4},\underline{\quad})$}}}&}&
a;
\\
\mathfrak{R}_{1}\colon b&
\xymatrix@C=40pt{\ar[r]^-{\text{\Small{$(\mathfrak{p}_{3},\underline{\quad})$}}}&}&
\delta(x,y)&
\xymatrix@C=40pt{\ar[r]^-{\text{\Small{$(\mathfrak{q}_{1},\delta(x,\underline{\;\;\;})))$}}}&}&
\delta(x,x)&
\xymatrix@C=40pt{\ar[r]^-{\text{\Small{$(\mathfrak{q}_{0},\delta(\underline{\;\;\;},x)))$}}}&}&
\delta(z,x);
\\
\mathfrak{R}_{2}\colon  a. & \textcolor{white}{\mathcal{A}^{\mathcal{A}^{\mathcal{A}^{\mathcal{A}^{\mathcal{A}^{\mathcal{A}}}}}}_{\mathcal{A}_{\mathcal{A}_{\mathcal{A}_{\mathcal{A}_{\mathcal{A}}}}}}}
\end{array}
\]
Let us note that (1) $\mathfrak{P}^{0,3}$ is a $(0110)$-path in $\boldsymbol{\mathcal{A}}$ of sort $0$, (2) $\mathfrak{R}_{0}$ is a $0$-path in $\boldsymbol{\mathcal{A}}$ of sort $0$, (3) $\mathfrak{R}_{1}$ is a $(011)$-path  in $\boldsymbol{\mathcal{A}}$ of sort $0$, and (4) $\mathfrak{R}_{2}$ is a $\lambda$-path in $\boldsymbol{\mathcal{A}}$ of sort $0$. In particular, $\mathfrak{R}_{0}$ is an echelon, $\mathfrak{R}_{1}$ is the composition of an echelon and an echelonless path, and $\mathfrak{R}_{1}$ is the $(1,0)$-identity path on $a$.

Consider the $(0011)$-path in $\boldsymbol{\mathcal{A}}$ of sort $0$ given by $\beta^{\mathbf{Pth}_{\boldsymbol{\mathcal{A}}}}((\mathfrak{R}_{j})_{j\in 3})$, i.e., the result of applying $\beta$  to the family of paths $(\mathfrak{R}_{j})_{j\in 3}$ in the $2$-sorted $\Sigma$-algebra $\mathbf{Pth}^{(0,1)}_{\boldsymbol{\mathcal{A}}}$. By Proposition~\ref{PPthAlg}, the path $\beta^{\mathbf{Pth}_{\boldsymbol{\mathcal{A}}}}((\mathfrak{R}_{j})_{j\in 3})$ is given by the following sequence
\[
\allowdisplaybreaks
\begin{array}{rcl}
\beta^{\mathbf{Pth}_{\boldsymbol{\mathcal{A}}}}((\mathfrak{R}_{j})_{j\in 3})\colon
\beta(\gamma(z, a),b,a)&
\xymatrix@C=100pt{\ar[r]^-{\text{\Small{$\left(\mathfrak{p}_{4}, 
\beta(\underline{\quad},b,a)
\right)$}}}&}&
\beta(a,b,a)
\\&
\xymatrix@C=100pt{\ar[r]^-{\text{\Small{$\left(\mathfrak{p}_{3}, 
\beta(a,\underline{\quad},a)
\right)$}}}&}&
\beta(a,\delta(x,y),a)
\\&
\xymatrix@C=100pt{\ar[r]^-{\text{\Small{$\left(\mathfrak{q}_{1}, 
\beta(a,\delta(x,\underline{\quad}),a)
\right)$}}}&}&
\beta(a,\delta(x,x),a)
\\&
\xymatrix@C=100pt{\ar[r]^-{\text{\Small{$\left(\mathfrak{q}_{0}, 
\beta(a,\delta(\underline{\quad},x),a)
\right)$}}}&}&
\beta(a,\delta(z,x),a).
\end{array}
\]

The paths that we can extract from $\beta^{\mathbf{Pth}_{\boldsymbol{\mathcal{A}}}}((\mathfrak{R}_{j})_{j\in 3})$ are, by Proposition~\ref{PRecov}, again, the components of the family of paths $(\mathfrak{R}_{j})_{j\in 3}$. However, the appearance is no longer chaotic. We see how the path extraction algorithm, from Lemma~\ref{LPthExtract}, applies to $\beta^{\mathbf{Pth}_{\boldsymbol{\mathcal{A}}}}((\mathfrak{R}_{j})_{j\in 3})$.
\begin{center}
\begin{tikzpicture}[
fletxa/.style={thick, ->}]
\tikzstyle{every state}=[very thick, draw=blue!50,fill=blue!20, inner sep=2pt,minimum size=20pt]
\node[] (0P) at (0,0) [] {$P_{8}$};
\node[] (0=) at (.7,0) [] {$=$};
\node[] (0s) at (1.2,0) [] {$\beta($};
\node[] (00) at (2.2,0) [] {$\color{white}{P^{i}_{\bb{\mathbf{s}}-1}}$};
\node[] () at (2.2,0) [] {$\gamma(z, a)$};
\node[] (0c0) at (3.4,0) [] {$,$};
\node[] (0j) at (4.5,0) [] {$\color{white}{P^{i}_{\bb{\mathbf{s}}-1}}$};
\node[] () at (4.5,0) [] {$b$};
\node[] (0c1) at (5.6,0) [] {$,$};
\node[] (0f) at (6.7,0) [] {$\color{white}{P^{i}_{\bb{\mathbf{s}}-1}}$};
\node[] () at (6.7,0) [] {$a$};
\node[] (0c0) at (7.6,0) [] {$)$};

\node[] (E0label) at (-1,-.75) [] {$\scriptstyle (\mathfrak{p}_{4},\beta(\underline{\quad},b,a))$};
\node[] (E0) at (0,-.75) [] {\rotatebox{-90}{$\mor$}};
\node[] (E00) at (2.2,-.75) [] {\rotatebox{-90}{$\mor$}};
\node[] (E00label) at (2.85,-.75) [] {$\scriptstyle (\mathfrak{p}_{4},\underline{\quad})$};
\node[] (E0j) at (4.5,-.75) [] {\rotatebox{-90}{$\,=$}};
\node[] (E0f) at (6.7,-.75) [] {\rotatebox{-90}{$\,=$}};

\node[] (1P) at (0,-1.5) [] {$Q_{9}$};
\node[] (1=) at (.7,-1.5) [] {$=$};
\node[] (1s) at (1.2,-1.5) [] {$\beta($};
\node[] (10) at (2.2,-1.5) [] {$\color{white}P^{i}_{\bb{\mathbf{s}}-1}$};
\node[] () at (2.2,-1.5) [] {$a$};
\node[] (1c0) at (3.4,-1.5) [] {$, $};
\node[] (1j) at (4.5,-1.5) [] {$\color{white}P^{i}_{\bb{\mathbf{s}}-1}$};
\node[] () at (4.5,-1.5) [] {$b$};
\node[] (1c1) at (5.6,-1.5) [] {$, $};
\node[] (1f) at (6.7,-1.5) [] {$\color{white}P^{i}_{\bb{\mathbf{s}}-1}$};
\node[] () at (6.7,-1.5) [] {$a$};
\node[] (1c0) at (7.6,-1.5) [] {$)$};

\node[] (E1label) at (-1,-2.25) [] {$\scriptstyle (\mathfrak{p}_{3},\beta(a,\underline{\quad},a))$};
\node[] (E1) at (0,-2.25) [] {\rotatebox{-90}{$\mor$}};
\node[] (E10) at (2.2,-2.25) [] {\rotatebox{-90}{$\,=$}};
\node[] (E1j) at (4.5,-2.25) [] {\rotatebox{-90}{$\mor$}};
\node[] (E1jlabel) at (5.15,-2.25) [] {$\scriptstyle (\mathfrak{p}_{3},\underline{\quad})$};
\node[] (E1f) at (6.7,-2.25) [] {\rotatebox{-90}{$\,=$}};

\node[] (2P) at (0,-3) [] {$Q_{10}$};
\node[] (2=) at (.7,-3) [] {$=$};
\node[] (2s) at (1.2,-3) [] {$\beta($};
\node[] (20) at (2.2,-3) [] {$\color{white}P^{i}_{\bb{\mathbf{s}}-1}$};
\node[] () at (2.2,-3) [] {$a$};
\node[] (2c0) at (3.4,-3) [] {$, $};
\node[] (2j) at (4.5,-3) [] {$\color{white}P^{i}_{\bb{\mathbf{s}}-1}$};
\node[] () at (4.5,-3) [] {$\delta(x,y)$};
\node[] (2c1) at (5.6,-3) [] {$, $};
\node[] (2f) at (6.7,-3) [] {$\color{white}P^{i}_{\bb{\mathbf{s}}-1}$};
\node[] () at (6.7,-3) [] {$a$};
\node[] (2c0) at (7.6,-3) [] {$)$};

\node[] (E2label) at (-1.3,-3.75) [] {$\scriptstyle (\mathfrak{q}_{1}, \beta(a,\delta(x,\underline{\quad}),a) )$};
\node[] (E2) at (0,-3.75) [] {\rotatebox{-90}{$\mor$}};
\node[] (E20) at (2.2,-3.75) [] {\rotatebox{-90}{$\,=$}};
\node[] (E2j) at (4.5,-3.75) [] {\rotatebox{-90}{$\mor$}};
\node[] (E2jlabel) at (5.35,-3.75) [] {$\scriptstyle (\mathfrak{q}_{1}, \delta(x,\underline{\;\;\;}))$};
\node[] (E2f) at (6.7,-3.75) [] {\rotatebox{-90}{$\,=$}};

\node[] (3P) at (0,-4.5) [] {$Q_{11}$};
\node[] (3=) at (.7,-4.5) [] {$=$};
\node[] (3s) at (1.2,-4.5) [] {$\beta($};
\node[] (30) at (2.2,-4.5) [] {$\color{white}P^{i}_{\bb{\mathbf{s}}-1}$};
\node[] () at (2.2,-4.5) [] {$ a$};
\node[] (3c0) at (3.4,-4.5) [] {$, $};
\node[] (3j) at (4.5,-4.5) [] {$\color{white}P^{i}_{\bb{\mathbf{s}}-1}$};
\node[] () at (4.5,-4.5) [] {$\delta(x,x)$};
\node[] (3c1) at (5.6,-4.5) [] {$, $};
\node[] (3f) at (6.7,-4.5) [] {$\color{white}P^{i}_{\bb{\mathbf{s}}-1}$};
\node[] () at (6.7,-4.5) [] {$a$};
\node[] (3c0) at (7.6,-4.5) [] {$)$};

\node[] (E3label) at (-1.3,-5.25) [] {$\scriptstyle (\mathfrak{q}_{0}, \beta(a,\delta(\underline{\quad},x),a) )$};
\node[] (E3) at (0,-5.25) [] {\rotatebox{-90}{$\mor$}};
\node[] (E30) at (2.2,-5.25) [] {\rotatebox{-90}{$\,=$}};
\node[] (E3jlabel) at (5.35,-5.25) [] {$\scriptstyle (\mathfrak{q}_{0}, \delta(\underline{\;\;\;},x))$};
\node[] (E3j) at (4.5,-5.25) [] {\rotatebox{-90}{$\mor$}};
\node[] (E3f) at (6.7,-5.25) [] {\rotatebox{-90}{$\,=$}};

\node[] (4P) at (0,-6) [] {$P_{12}$};
\node[] (4=) at (.7,-6) [] {$=$};
\node[] (4s) at (1.2,-6) [] {$\beta($};
\node[] (40) at (2.2,-6) [] {$\color{white}P^{i}_{\bb{\mathbf{s}}-1}$};
\node[] () at (2.2,-6) [] {$ a$};
\node[] (4c0) at (3.4,-6) [] {$, $};
\node[] (4j) at (4.5,-6) [] {$\color{white}P^{i}_{\bb{\mathbf{s}}-1}$};
\node[] () at (4.5,-6) [] {$\delta(z,x)$};
\node[] (4c1) at (5.6,-6) [] {$, $};
\node[] (4f) at (6.7,-6) [] {$\color{white}P^{i}_{\bb{\mathbf{s}}-1}$};
\node[] () at (6.7,-6) [] {$a$};
\node[] (4c0) at (7.6,-6) [] {$)$};

\tikzset{encercla/.style={draw=black, line width=.5pt, inner sep=0pt, rectangle, rounded corners}};
\node [encercla,  fit=(00)  ] {} ;
\node [encercla,  fit=(10) (20) (30) (40)  ] {} ;
\node [encercla,  fit=(0j) (1j)  ] {} ;
\node [encercla,  fit= (2j)  ] {} ;
\node [encercla,  fit= (3j)  ] {} ;
\node [encercla,  fit= (4j)  ] {} ;
\node [encercla,  fit=(0f) (1f) (2f) (3f) (4f)  ] {} ;
\end{tikzpicture}
\end{center}
The path $\beta^{\mathbf{Pth}_{\boldsymbol{\mathcal{A}}}}((\mathfrak{R}_{j})_{j\in 3})$ does not follow a leftmost innermost derivation strategy. This is because the component $\mathfrak{R}_{1}$ is a path that does not satisfy this condition. Although the intermediate terms may have been changed, this path has the same source and target as $\mathfrak{P}^{7,11}$.
\end{example}
\chapter{
\texorpdfstring
{An Artinian order on $\coprod\mathrm{Pth}_{\boldsymbol{\mathcal{A}}}$}
{An Artinian order on paths}
}\label{S1C}


In this chapter we define on $\coprod\mathrm{Pth}_{\boldsymbol{\mathcal{A}}}$, the coproduct of $\mathrm{Pth}_{\boldsymbol{\mathcal{A}}}$---formed by all labelled paths $(\mathfrak{P},s)$, with $s\in S$ and $\mathfrak{P}\in \mathrm{Pth}_{\boldsymbol{\mathcal{A}},s}$---, an Artinian order (see Definition~\ref{DPosArt}), which will allow us to justify both proofs by Artinian induction and definitions by Artinian recursion, which we will use in subsequent chapters. Moreover, we specify the minimals of the corresponding Artinian ordered set, prove that $\mathrm{ip}^{(1,0)\sharp}[\mathrm{T}_{\Sigma}(X)]$ is a lower set and that several mappings from and to $\coprod\mathrm{Pth}_{\boldsymbol{\mathcal{A}}}$ are order-preserving, order-reflecting or order embeddings.

\begin{restatable}{definition}{DOrd}
\label{DOrd}
\index{partial order!first-order!$\leq_{\mathbf{Pth}_{\boldsymbol{\mathcal{A}}}}$}
We let $\prec_{\mathbf{Pth}_{\boldsymbol{\mathcal{A}}}}$ denote the binary relation on $\coprod \mathrm{Pth}_{\boldsymbol{\mathcal{A}}}$ consisting of the ordered pairs $((\mathfrak{Q},t),(\mathfrak{P},s))\in(\coprod \mathrm{Pth}_{\boldsymbol{\mathcal{A}}})^{2}$ for which one of the following conditions holds
\begin{enumerate}
\item $\mathfrak{P}$ and $\mathfrak{Q}$ are $(1,0)$-identity paths of the form
\begin{align*}
\mathfrak{P}&=\mathrm{ip}^{(1,0)\sharp}_{s}\left(P\right),
&
\mathfrak{Q}&=\mathrm{ip}^{(1,0)\sharp}_{t}\left(Q\right),
\end{align*}
for some terms $P\in\mathrm{T}_{\Sigma}(X)_{s}$ and $Q\in\mathrm{T}_{\Sigma}(X)_{t}$ and the  inequality $$
\left(Q,t
\right)
<_{\mathbf{T}_{\Sigma}(X)}
\left(P,s
\right)
$$ holds, where $\leq_{\mathbf{T}_{\Sigma}(X)}$ is the subterm preorder on $\coprod\mathrm{T}_{\Sigma}(X)$, see Remark~\ref{RTermOrd}.
\item $\mathfrak{P}$ is a path of length strictly greater than one containing at least one echelon, and if its first echelon occurs at position $i\in\bb{\mathfrak{P}}$, then
\subitem(2.1) if $i=0$, then $\mathfrak{Q}$ is equal to $\mathfrak{P}^{0,0}$ or $\mathfrak{P}^{1,\bb{\mathfrak{P}}-1}$,
\subitem(2.2) if $i>0$, then $\mathfrak{Q}$ is equal to $\mathfrak{P}^{0,i-1}$ or $\mathfrak{P}^{i,\bb{\mathfrak{P}}-1}$;
\item $\mathfrak{P}$ is an echelonless path and $\mathfrak{Q}$ is one of the paths extracted from $\mathfrak{P}$ in virtue of Lemma~\ref{LPthExtract}.
\end{enumerate}
We will denote by $\leq_{\mathbf{Pth}_{\boldsymbol{\mathcal{A}}}}$ the reflexive and transitive closure of $\prec_{\mathbf{Pth}_{\boldsymbol{\mathcal{A}}}}$, i.e., the preorder on $\coprod \mathrm{Pth}_{\boldsymbol{\mathcal{A}}}$  generated by $\prec_{\mathbf{Pth}_{\boldsymbol{\mathcal{A}}}}$, and by $<_{\mathbf{Pth}_{\boldsymbol{\mathcal{A}}}}$ the transitive closure of $\prec_{\mathbf{Pth}_{\boldsymbol{\mathcal{A}}}}$.
\end{restatable}

\begin{figure}
\begin{tikzpicture}[scale=1.1, baseline=0cm]
\tikzstyle{end} =  [regular polygon,regular polygon sides=5, draw, shape border rotate=180,  minimum width=.8cm, minimum height=.8cm, text centered, draw=black, fill=purple!20!white]
\tikzstyle{start} =  [regular polygon,regular polygon sides=5, draw, minimum width=.8cm, minimum height=.8cm, text centered, draw=black, fill=purple!20!white]
\tikzstyle{true} = [rectangle, rounded corners, minimum width=.7cm, minimum height=.7cm, text centered, draw=black, fill=green!20!white]
\tikzstyle{false} = [rectangle, rounded corners, minimum width=.7cm, minimum height=.7cm, text centered, draw=black, fill=red!20!white]
\tikzstyle{io} = [trapezium, trapezium left angle=70, trapezium right angle=110, minimum width=.55cm, minimum height=.55cm, text centered, draw=black, fill=purple!20!white]
\tikzstyle{process} = [rectangle, minimum width=.7cm, minimum height=.7cm, text centered, draw=black, fill=yellow!20!white]
\tikzstyle{decision} = [diamond, minimum width=.8cm, minimum height=.8cm, text centered, draw=black, fill=blue!20!white]
\tikzstyle{arrow} = [thick,->,>=stealth]

\node (A1) at (0.5,-1) [start] {$\scriptstyle 1$};
\node (A2) at (0.5,-2) [decision] {$\scriptstyle 2$};
\node (A3) at (-1,-3) [decision] {$\scriptstyle 3$};
\node (A4) at (2,-3) [decision] {$\scriptstyle 4$};
\node (A5) at (-2,-4) [process] {$\scriptstyle 5$};
\node (F0) at (-1,-4) [false] {$\scriptstyle \mathrm{F}$};
\node (A6) at (1,-4) [decision] {$\scriptstyle 6$};
\node (F1) at (1,-5) [false] {$\scriptstyle \mathrm{F}$};
\node (A7) at (3,-4) [decision] {$\scriptstyle 7$};
\node (A8) at (-2,-5) [end] {$\scriptstyle 8$};
\node (A9) at (0,-5) [process] {$\scriptstyle 9$};
\node (A10) at (0,-6) [decision] {$\scriptstyle 10$};
\node (A11) at (-1,-7) [decision] {$\scriptstyle 11$};
\node (V1) at (-3,-6) [true] {$\scriptstyle \mathrm{T}$};
\node (F2) at (-2,-6) [false] {$\scriptstyle \mathrm{F}$};
\node (A12) at (1,-7) [decision] {$\scriptstyle 12$};
\node (V2) at (-2,-8) [true] {$\scriptstyle \mathrm{T}$};
\node (F3) at (-1,-8) [false] {$\scriptstyle \mathrm{F}$};
\node (V3) at (0,-8) [true] {$\scriptstyle \mathrm{T}$};
\node (F4) at (1,-8) [false] {$\scriptstyle \mathrm{F}$};
\node (V4) at (2,-5) [true] {$\scriptstyle \mathrm{T}$};
\node (F5) at (3,-5) [false] {$\scriptstyle \mathrm{F}$};

\draw [arrow] (A1) -- (A2);
\draw [arrow] (A2) -| node[anchor=south west] {yes} (A3);
\draw [arrow] (A2) -| node[anchor=south east] {no} (A4);
\draw [arrow] (A3) -| node[anchor=south west] {yes} (A5);
\draw [arrow] (A3) -- node[anchor=west] {no} (F0);
\draw [arrow] (A4) -| node[anchor=south west] {yes} (A6);
\draw [arrow] (A4) -| node[anchor=south east] {no} (A7);
\draw [arrow] (A6) -| node[anchor=south west] {yes} (A9);
\draw [arrow] (A6) -- node[anchor=west] {no} (F1);
\draw [arrow] (A7) -| node[anchor=south west] {yes} (V4);
\draw [arrow] (A7) -- node[anchor=west] {no} (F5);
\draw [arrow] (A5) -- (A8);
\draw [arrow] (A9) -- (A10);
\draw [arrow] (A8) -| node[anchor=south west] {yes} (V1);
\draw [arrow] (A8) -- node[anchor=west] {no} (F2);
\draw [arrow] (A10) -| node[anchor=south west] {yes} (A11);
\draw [arrow] (A10) -| node[anchor=south east] {no} (A12);
\draw [arrow] (A11) -| node[anchor=south west] {yes} (V2);
\draw [arrow] (A11) -- node[anchor=west] {no} (F3);
\draw [arrow] (A12) -| node[anchor=south west] {yes} (V3);
\draw [arrow] (A12) -- node[anchor=west] {no} (F4);

\end{tikzpicture}
\begin{center}
\begin{tabu} to 1.2\textwidth {  r  X[l]  }
\\
1.&Let $(\mathfrak{Q},t),(\mathfrak{P},s)$ be paths in $\coprod\mathrm{Pth}_{\boldsymbol{\mathcal{A}}}$\\
2.&Is $\mathfrak{P}$ a $(1,0)$-identity path?\\
3.&Is $\mathfrak{Q}$ a $(1,0)$-identity path?\\
4.&Does $\mathfrak{P}$ contain at least one echelon?\\
5.&Let $\mathfrak{Q}=\mathrm{ip}^{(1,0)\sharp}_{t}(Q)$ and
$\mathfrak{P}=\mathrm{ip}^{(1,0)\sharp}_{s}(P)$\\
6.&Is $\bb{\mathfrak{P}}>1$?\\
7.&Is $\mathfrak{Q}$ one of the paths we can extract from $\mathfrak{P}$?\\
8.&Is $(Q,t)\leq_{\mathbf{T}_{\Sigma}(X)}(P,s)$?\\
9.&Let $i\in\bb{\mathfrak{P}}$ be the first index for which $\mathfrak{P}^{i,i}$ is an echelon\\
10.&Is $i=0$?\\
11.&Is $\mathfrak{Q}$ equal to $\mathfrak{P}^{0,0}$ or $\mathfrak{P}^{1,\bb{\mathfrak{P}}-1}$?\\
12.&Is $\mathfrak{Q}$ equal to $\mathfrak{P}^{0,i-1}$ or $\mathfrak{P}^{i,\bb{\mathfrak{P}}-1}$?\\
\\
\end{tabu}
\end{center}
\begin{tikzpicture}[scale=1.1, baseline=0cm]
\tikzstyle{true} = [rectangle, rounded corners, minimum width=.7cm, minimum height=.7cm, text centered, draw=black, fill=green!20!white]
\tikzstyle{false} = [rectangle, rounded corners, minimum width=.7cm, minimum height=.7cm, text centered, draw=black, fill=red!20!white]

\node () at (0,0) [true] {$\scriptstyle \mathrm{T}$};
\node () at (2,0) [] {$~(\mathfrak{Q},t)\prec_{\mathbf{Pth}_{\boldsymbol{\mathcal{A}}}}(\mathfrak{P},s)$};
\node () at (5,0) [false] {$\scriptstyle \mathrm{F}$};
\node () at (7,0) [] {$~(\mathfrak{Q},t)\not\prec_{\mathbf{Pth}_{\boldsymbol{\mathcal{A}}}}(\mathfrak{P},s)$};
\end{tikzpicture}
\caption{Decision flowchart for $(\mathfrak{Q},t)\prec_{\mathbf{Pth}_{\boldsymbol{\mathcal{A}}}}(\mathfrak{P},s)$.}\label{FFlow}
\end{figure}

\begin{restatable}{remark}{RIrrefl}
\label{RIrrefl}
The binary relation $\prec_{\mathbf{Pth}_{\boldsymbol{\mathcal{A}}}}$ on $\coprod \mathrm{Pth}_{\boldsymbol{\mathcal{A}}}$ is irreflexive.
\end{restatable}

\begin{restatable}{remark}{ROrd}
\label{ROrd}
The preorder $\leq_{\mathbf{Pth}_{\boldsymbol{\mathcal{A}}}}$ on $\coprod \mathrm{Pth}_{\boldsymbol{\mathcal{A}}}$ is $\bigcup_{n\in\mathbb{N}}\prec_{\mathbf{Pth}_{\boldsymbol{\mathcal{A}}}}^{n}$, where $\prec_{\mathbf{Pth}_{\boldsymbol{\mathcal{A}}}}^{0}$ is the diagonal of $\coprod \mathrm{Pth}_{\boldsymbol{\mathcal{A}}}$ and, for $n\in\mathbb{N}$, 
$$\prec_{\mathbf{Pth}_{\boldsymbol{\mathcal{A}}}}^{n+1} = \prec_{\mathbf{Pth}_{\boldsymbol{\mathcal{A}}}}^{n}\circ \prec_{\mathbf{Pth}_{\boldsymbol{\mathcal{A}}}}.$$ 

Thus, for every $((\mathfrak{Q},t),(\mathfrak{P},s))\in (\coprod \mathrm{Pth}_{\boldsymbol{\mathcal{A}}})^{2}$, $(\mathfrak{Q},t)\leq_{\mathbf{Pth}_{\boldsymbol{\mathcal{A}}}}(\mathfrak{P},s)$ if and only if $s = t$ and $\mathfrak{Q} = \mathfrak{P}$ or there exists  a natural number $m\in\mathbb{N}-\{0\}$, a word $\mathbf{w}\in S^{\star}$ of length $\bb{\mathbf{w}}=m+1$, and a family of paths $(\mathfrak{R}_{k})_{k\in\bb{\mathbf{w}}}$ in $\mathrm{Pth}_{\boldsymbol{\mathcal{A}},{\mathbf{w}}}$, such that $w_{0}=t$,  $\mathfrak{R}_{0}=\mathfrak{Q}$, $w_{m}=s$, $\mathfrak{R}_{m}=\mathfrak{P}$ and, for every $k\in m$, $(\mathfrak{R}_{k},w_{k})\prec_{\mathbf{Pth}_{\boldsymbol{\mathcal{A}}}}(\mathfrak{R}_{k+1},w_{k+1})$. Moreover, $<_{\mathbf{Pth}_{\boldsymbol{\mathcal{A}}}}$, the transitive closure of $\prec_{\mathbf{Pth}_{\boldsymbol{\mathcal{A}}}}$, is $\bigcup_{n\in\mathbb{N}-
\{0\}}\prec_{\mathbf{Pth}_{\boldsymbol{\mathcal{A}}}}^{n}$.
\end{restatable}

In the following proposition we show that $\mathrm{Min}(\coprod \mathrm{Pth}_{\boldsymbol{\mathcal{A}}}, \leq_{\mathbf{Pth}_{\boldsymbol{\mathcal{A}}}})$, the set of the minimals of the preordered set $(\coprod \mathrm{Pth}_{\boldsymbol{\mathcal{A}}}, \leq_{\mathbf{Pth}_{\boldsymbol{\mathcal{A}}}})$, consists of the $(1,0)$-identity paths on minimal terms, i.e., on terms in $\mathrm{B}^{0}_{\Sigma}(X)$, together with the echelon paths in $\mathrm{Pth}_{\boldsymbol{\mathcal{A}},s}$, labelled by their corresponding sorts in $S$.

\begin{restatable}{proposition}{PMinimal}
\label{PMinimal}
For the preordered set $(\coprod \mathrm{Pth}_{\boldsymbol{\mathcal{A}}}, \leq_{\mathbf{Pth}_{\boldsymbol{\mathcal{A}}}})$ we have that
\begin{multline*}
\textstyle
\mathrm{Min}
\left(
\coprod \mathrm{Pth}_{\boldsymbol{\mathcal{A}}}, \leq_{\mathbf{Pth}_{\boldsymbol{\mathcal{A}}}}
\right) 
\\
\textstyle
= 
\left(\coprod\mathrm{ip}^{(1,0)\sharp}\left[
\mathrm{Min}
\left(
\coprod \mathrm{T}_{\Sigma}(X), \leq_{\mathbf{T}_{\Sigma}(X)}
\right)
\right]\right)
\cup
\left(\coprod\mathrm{Ech}^{(1,\mathcal{A})}\left[\mathcal{A}
\right]\right).
\end{multline*}
\end{restatable}

\begin{proof}
We first prove that both labelled $(1,0)$-identity paths on minimal terms, i.e., terms in $\mathrm{B}^{0}_{\Sigma}(X)$, and  labelled echelon paths of any sort are minimal elements.

Let $s$ be a sort in $S$ and $P$ a term in $\mathrm{T}_{\Sigma}(X)_{s}$. Assume that $(P,s)$ is a minimal element in $(\coprod \mathrm{T}_{\Sigma}(X), \leq_{\mathbf{T}_{\Sigma}(X)})$.   Consider the $(1,0)$-identity path on $P$, i.e., $\mathrm{ip}^{(1,0)\sharp}_{s}(P)$. Let us note that $\mathrm{ip}^{(1,0)\sharp}_{s}(P)$ is a path of length $0$. Therefore, for every sort $t\in S$ and every path $\mathfrak{Q}\in\mathrm{Pth}_{\boldsymbol{\mathcal{A}},t}$, we have that $(\mathfrak{Q},t)\prec_{\mathbf{Pth}_{\boldsymbol{\mathcal{A}}}} (\mathrm{ip}^{(1,0)\sharp}_{s}(P),s)$ if, and only if, $\mathfrak{Q}$ is a $(1,0)$-identity path of the form $\mathfrak{Q}=\mathrm{ip}^{(1,0)}_{t}(Q)$, for some term $Q\in\mathrm{T}_{\Sigma}(X)_{t}$ and 
$$
\left(
Q ,t\right)
\leq_{\mathbf{T}_{\Sigma}(X)}
\left(P,s
\right).
$$

Since $(P,s)$ is a minimal element in $(\coprod\mathrm{T}_{\Sigma}(X),\leq_{\mathbf{T}_{\Sigma}(X)})$, we have that $s=t$ and $Q=P$. Therefore, we have that 
$$
\mathfrak{Q}=
\mathrm{ip}^{(1,0)\sharp}_{t}\left(
Q
\right)
=
\mathrm{ip}^{(1,0)\sharp}_{s}\left(
P
\right)
=
\mathfrak{P}.
$$
It follows that $(\mathrm{ip}^{(1,0)\sharp}_{s}(P),s)$ is minimal in $(\coprod \mathrm{Pth}_{\boldsymbol{\mathcal{A}}}, \leq_{\mathbf{Pth}_{\boldsymbol{\mathcal{A}}}})$.

Let $s$ be a sort in $S$ and $\mathfrak{p}$ a rewrite rule in $\mathcal{A}_{s}$ of the form $\mathfrak{p}=(M,N)$. Consider the echelon associated to $\mathfrak{p}$
$$
\xymatrix@C=90pt{
\mathrm{ech}^{(1,\mathcal{A})}_{s}(\mathfrak{p})\colon
M
\ar@{->}[r]^-{\text{\Small{($\mathfrak{p}$, $\mathrm{id}^{\mathrm{T}_{\Sigma}(X)_{s}}$)}}}
&
N
}
.
$$

Let us note that 
$\mathrm{ech}^{(1,\mathcal{A})}_{s}(\mathfrak{p})$ is neither a $(1,0)$-identity path, nor a  path of length strictly greater than one  nor an echelonless path. Therefore, for every sort $t\in S$ and every path $\mathfrak{Q}\in\mathrm{Pth}_{\boldsymbol{\mathcal{A}},t}$, we have that $(\mathfrak{Q},t)\nprec_{\mathbf{Pth}_{\boldsymbol{\mathcal{A}}}} (\mathrm{ech}^{(1,\mathcal{A})}_{s}(\mathfrak{p}),s)$. Consequently, no element of $\coprod \mathrm{Pth}_{\boldsymbol{\mathcal{A}}}$ different from $(\mathrm{ech}^{(1,\mathcal{A})}_{s}(\mathfrak{p}),s)$ exists which $<_{\mathbf{Pth}_{\boldsymbol{\mathcal{A}}}}$-precedes $(\mathrm{ech}^{(1,\mathcal{A})}_{s}(\mathfrak{p}),s)$. From this it follows that $(\mathrm{ech}^{(1,\mathcal{A})}_{s}(\mathfrak{p}),s)$ is minimal in $(\coprod \mathrm{Pth}_{\boldsymbol{\mathcal{A}}}, \leq_{\mathbf{Pth}_{\boldsymbol{\mathcal{A}}}})$.

Conversely, we prove that every minimal element of $(\coprod \mathrm{Pth}_{\boldsymbol{\mathcal{A}}}, \leq_{\mathbf{Pth}_{\boldsymbol{\mathcal{A}}}})$ is a labelled $(1,0)$-identity path on a minimal term or a labelled echelon.

Let $s$ be a sort in $S$ and $\mathfrak{P}$ a path in $\mathrm{Pth}_{\boldsymbol{\mathcal{A}},s}$. Assume that $(\mathfrak{P},s)$ is a minimal element of $(\coprod \mathrm{Pth}_{\boldsymbol{\mathcal{A}}}, \leq_{\mathbf{Pth}_{\boldsymbol{\mathcal{A}}}})$. Then we want to show that $(\mathfrak{P},s)$ belongs to $\coprod\mathrm{Ech}[\mathcal{A}]$ or to
$(\coprod\mathrm{ip}^{(1,0)\sharp}[
\mathrm{Min}
(
\coprod \mathrm{T}_{\Sigma}(X), \leq_{\mathbf{T}_{\Sigma}(X)}
)
])$.

We consider three cases. Either (1) $\mathfrak{P}$ is a $(1,0)$-identity path or (2) $\mathfrak{P}$ is a path of length at least one containing echelons or (3) $\mathfrak{P}$ is an echelonless path.

If (1) $\mathfrak{P}$ is a $(1,0)$-identity path of the form $\mathfrak{P}=\mathrm{ip}^{(1,0)\sharp}_{s}(P)$, for a unique  term  $P$ in $\mathrm{T}_{\Sigma}(X)_{s}$, then either (1.1)  $(P,s)$ is a minimal pair or (1.2) $(P,s)$ is not a minimal pair in the partially ordered set  $(\coprod \mathrm{T}_{\Sigma}(X), \leq_{\mathbf{T}_{\Sigma}(X)})$.

If (1.1) $(P,s)$ is a minimal pair in  $(\coprod \mathrm{T}_{\Sigma}(X), \leq_{\mathbf{T}_{\Sigma}(X)})$, then the statement follows because $(\mathfrak{P},s)$ is a labelled pair in 
$$
\textstyle
\left(\coprod\mathrm{ip}^{(1,0)\sharp}\left[
\mathrm{Min}
\left(
\coprod \mathrm{T}_{\Sigma}(X), \leq_{\mathbf{T}_{\Sigma}(X)}
\right)
\right]\right).
$$

If (1.2) $(P,s)$ is not a minimal pair in  $(\coprod \mathrm{T}_{\Sigma}(X), \leq_{\mathbf{T}_{\Sigma}(X)})$, then one can find a sort $t\in S$ and a term $Q$ in $\mathrm{T}_{\Sigma}(X)_{t}$ for which
$$
\left(Q, t
\right)
<_{\mathbf{T}_{\Sigma}(X)}
\left(P, s
\right).
$$

Consider the path $\mathfrak{Q}$ in $\mathrm{Pth}_{\boldsymbol{\mathcal{A}},t}$ given by $\mathfrak{Q}=\mathrm{ip}^{(1,0)\sharp}_{t}(Q)$. Then, taking into account the above strict inequality and according to Definition~\ref{DOrd}, we have that 
$$
\left(\mathfrak{Q},t
\right)
<_{\mathbf{Pth}_{\boldsymbol{\mathcal{A}}}}
\left(\mathfrak{P},s
\right),
$$
contradicting that $(\mathfrak{P},s)$ is minimal.  This excludes case (1.2).

Now consider the case (2) in which  $\mathfrak{P}$ is a path of length at least one containing an echelon. Then either (2.1) the length of $\mathfrak{P}$ is one or (2.2) the length $\mathfrak{P}$ is strictly greater than one.

If (2.1), then $\mathfrak{P}$ is an echelon and the statement holds because $\mathfrak{P}$ is an element of $\mathrm{Ech}^{(1,\mathcal{A})}_{s}[\mathcal{A}_{s}]$.

If (2.2), then $\mathfrak{P}$ is a path of length strictly greater than one containing at least one echelon. Then there exists an index $i\in\bb{\mathfrak{P}}$ such that the one-step subpath $\mathfrak{P}^{i,i}$ is the first echelon appearing in $\mathfrak{P}$.  Regarding the nature of $i$, we have the following possibilities

If  $i=0$ then, according to Definition~\ref{DOrd}, we have that $(\mathfrak{P}^{0,0},s)$ $\prec_{\mathbf{Pth}_{\boldsymbol{\mathcal{A}}}}$-precedes $
(\mathfrak{P},s)$. Let us note that $\mathfrak{P}^{0,0}$ cannot be equal to $\mathfrak{P}$ because it has length $1$. This contradicts that $(\mathfrak{P},s)$ is minimal.

If  $i>0$, then, according to Definition~\ref{DOrd}, we have that $(\mathfrak{P}^{0,i-1},s)$ $\prec_{\mathbf{Pth}_{\boldsymbol{\mathcal{A}}}}$-precedes $
(\mathfrak{P},s)$. Let us note that $\mathfrak{P}^{0,i-1}$ cannot be equal to $\mathfrak{P}$ because it has smaller length. This contradicts that $(\mathfrak{P},s)$ is minimal.  

In any case, this excludes paths in case (2.2) of being minimal.

Let us consider, finally, the case (3) of $\mathfrak{P}$ being an echelonless path, then this path is suitable for the application of the path extraction algorithm. Let $\mathfrak{Q}$ be a path of sort $t\in S$ extracted from $\mathfrak{P}$, according to Lemma~\ref{LPthExtract}, then, by Definition~\ref{DOrd}, 
$(\mathfrak{Q},t)$ $\prec_{\mathbf{Pth}_{\boldsymbol{\mathcal{A}}}}$-precedes $
(\mathfrak{P},s)$. It follows from Lemma~\ref{LPthExtract} that $\mathfrak{Q}$ cannot be equal to $\mathfrak{P}$. 

This excludes paths in case (3) of being minimal.
\end{proof}

In the following proposition we prove that the coproduct of the $S$-sorted set of $(1,0)$-identity paths is a lower set of $(\coprod\mathrm{Pth}_{\boldsymbol{\mathcal{A}}}, \leq_{\mathbf{Pth}_{\boldsymbol{\mathcal{A}}}})$ and that this subset reflects the preorder arising from $\mathrm{T}_{\Sigma}(X)$.

\begin{restatable}{proposition}{PLower}
\label{PLower} The subset $\coprod\mathrm{ip}^{(1,0)\sharp}[\mathrm{T}_{\Sigma}(X)]$ of $\coprod\mathrm{Pth}_{\boldsymbol{\mathcal{A}}}$ is a lower set of the ordered set $(\coprod\mathrm{Pth}_{\boldsymbol{\mathcal{A}}}, \leq_{\mathbf{Pth}_{\boldsymbol{\mathcal{A}}}})$. That is, for every pair of sorts $s,t\in S$ and paths $\mathfrak{P}\in\mathrm{Pth}_{\boldsymbol{\mathcal{A}},s}$ and $\mathfrak{Q}\in\mathrm{Pth}_{\boldsymbol{\mathcal{A}},t}$, if
\begin{enumerate}
\item $(\mathfrak{Q},t)\leq_{\mathbf{Pth}_{\boldsymbol{\mathcal{A}}}}(\mathfrak{P},s)$ and
\item $\mathfrak{P}$ is a $(1,0)$-identity path, 
\end{enumerate}
then $\mathfrak{Q}$ is also a $(1,0)$-identity  path. Moreover, if $\mathfrak{P}$ and $\mathfrak{Q}$ have the form
\begin{align*}
\mathfrak{P}&=
\mathrm{ip}^{(1,0)\sharp}_{s}\left(
P
\right),
&
\mathfrak{Q}&=
\mathrm{ip}^{(1,0)\sharp}_{t}\left(
Q
\right),
\end{align*}
for some terms $P\in\mathrm{T}_{\Sigma}(X)_{s}$ and $Q\in\mathrm{T}_{\Sigma}(X)_{t}$, then the following inequality holds
$$
\left(Q,t\right)
\leq_{\mathbf{T}_{\Sigma}(X)}
\left(P,s\right).
$$
\end{restatable}
\begin{proof}
Let $(\mathfrak{Q},t)$, $(\mathfrak{P},s)$ be two pairs in $\coprod\mathrm{Pth}_{\boldsymbol{\mathcal{A}}}$ with $(\mathfrak{Q},t)\leq_{\mathbf{Pth}_{\boldsymbol{\mathcal{A}}}}(\mathfrak{P},s)$. Assume that $\mathfrak{P}$ is a $(1,0)$-identity path. 

Since $(\mathfrak{Q},t)\leq_{\mathbf{Pth}_{\boldsymbol{\mathcal{A}}}}(\mathfrak{P},s)$, we have, by Remark~\ref{ROrd}, that either $s=t$ and $\mathfrak{Q}=\mathfrak{P}$ or there exists a natural number $m\in\mathbb{N}-\{0\}$, a word $\mathbf{w}\in S^{\star}$ of length $\bb{\mathbf{w}}=m+1$, and a family of paths $(\mathfrak{R}_{k})_{k\in\bb{\mathbf{w}}}$ in $\mathrm{Pth}_{\boldsymbol{\mathcal{A}},\mathbf{w}}$ such that $w_{0}=t$, $\mathfrak{R}_{0}=\mathfrak{Q}$, $w_{m}=s$, $\mathfrak{R}_{m}=\mathfrak{P}$ and, for every $k\in m$, 
$
(\mathfrak{R}_{k}, w_{k})
\prec_{\mathbf{Pth}_{\boldsymbol{\mathcal{A}}}}
(\mathfrak{R}_{k+1}, w_{k+1}).
$

If $s=t$ and $\mathfrak{Q}=\mathfrak{P}$ the statement trivially holds. Therefore, it suffices to prove the proposition for a non-trivial sequence instantiating that 
$$\left(\mathfrak{Q},t
\right)\leq_{\mathbf{Pth}_{\boldsymbol{\mathcal{A}}}}
\left(\mathfrak{P},s
\right).$$

The proof is by induction on $m\in\mathbb{N}-\{0\}$.

\textsf{Base step of the induction.}

For $m=1$ we have that $(\mathfrak{Q},t)\prec_{\mathbf{Pth}_{\boldsymbol{\mathcal{A}}}}(\mathfrak{P},s)$. Then, according to Definition~\ref{DOrd}, since $\mathfrak{P}$ is a $(1,0)$-identity  path, we have that the unique possibility for $\mathfrak{Q}$ is that it is also a $(1,0)$-identity path. 

Moreover, if the  paths $\mathfrak{Q}$ and $\mathfrak{P}$ have the form
\begin{align*}
\mathfrak{P}&=
\mathrm{ip}^{(1,0)\sharp}_{s}\left(
P
\right),
&
\mathfrak{Q}&=
\mathrm{ip}^{(1,0)\sharp}_{t}
\left(
Q
\right),
\end{align*}
for some terms $P\in\mathrm{T}_{\Sigma}(X)_{s}$ and $Q\in\mathrm{T}_{\Sigma}(X)_{t}$, then the  inequality also holds
$
\left(Q,t
\right)
\leq_{\mathbf{T}_{\Sigma}(X)}
\left(P,s
\right).
$

\textsf{Inductive step of the induction.}

Let us suppose that the statement holds for sequences of length up to $m\in\mathbb{N}-\{0\}$. That is, if $\mathbf{w}$ is a word in $S^{\star}$ of length $\bb{\mathbf{w}}=m+1$ and $(\mathfrak{R}_{k})_{k\in\bb{\mathbf{w}}}$ is a family of paths in $\mathrm{Pth}_{\boldsymbol{\mathcal{A}},\mathbf{w}}$ such that $w_{0}=t$, $\mathfrak{R}_{0}=\mathfrak{Q}$, $w_{m}=s$, $\mathfrak{R}_{m}=\mathfrak{P}$ and, for every $k\in m$, $(\mathfrak{R}_{k},w_{k})\prec_{\mathbf{Pth}_{\boldsymbol{\mathcal{A}}}} (\mathfrak{R}_{k+1}, w_{k+1})$ then, if $\mathfrak{P}$ is a $(1,0)$-identity path, it follows that $\mathfrak{Q}$ is also a $(1,0)$-identity path. 

Moreover, if the paths $\mathfrak{Q}$ and $\mathfrak{P}$ have the form
\begin{align*}
\mathfrak{P}&=
\mathrm{ip}^{(1,0)\sharp}_{s}\left(
P
\right),
&
\mathfrak{Q}&=
\mathrm{ip}^{(1,0)\sharp}_{t}\left(
Q
\right),
\end{align*}
for some terms $P\in\mathrm{T}_{\Sigma}(X)_{s}$ and $Q\in\mathrm{T}_{\Sigma}(X)_{t}$, then the following inequality also holds
$
\left(Q,t\right)
\leq_{\mathbf{T}_{\Sigma}(X)}
\left(P,s\right).
$

Now, let $\mathbf{w}$ be a word in $S^{\star}$ of length $\bb{\mathbf{w}}=m+2$ and $(\mathfrak{R}_{k})_{k\in\bb{\mathbf{w}}}$ a family of paths in $\mathrm{Pth}_{\boldsymbol{\mathcal{A}},\mathbf{w}}$ such that $w_{0}=t$, $\mathfrak{R}_{0}=\mathfrak{Q}$, $w_{m+1}=s$, $\mathfrak{R}_{m+1}=\mathfrak{P}$ and, for every $k\in m+1$, $(\mathfrak{R}_{k},w_{k})\prec_{\mathbf{Pth}_{\boldsymbol{\mathcal{A}}}} (\mathfrak{R}_{k+1}, w_{k+1})$. 

Consider the final subsequence $(\mathfrak{R}_{k})_{k\in [1,m+1]}$. This is a sequence of length $m$ instantiating that $(\mathfrak{R}_{1},w_{1})\leq_{\mathbf{Pth}_{\boldsymbol{\mathcal{A}}}} (\mathfrak{P},s)$. Now, since $\mathfrak{P}$ is a $(1,0)$-identity path, it follows, by induction, that $\mathfrak{R}_{1}$ is also a $(1,0)$-identity  path. 

Moreover, if the  paths $\mathfrak{R}_{1}$ and $\mathfrak{P}$ have the form
\begin{align*}
\mathfrak{P}&=
\mathrm{ip}^{(1,0)\sharp}_{s}\left(
P
\right),
&
\mathfrak{R}_{1}&=
\mathrm{ip}^{(1,0)\sharp}_{w_{1}}\left(
R_{1}
\right),
\end{align*}
for some path terms $P\in\mathrm{T}_{\Sigma}(X)_{s}$ and $R_{1}\in\mathrm{T}_{\Sigma}(X)_{w_{1}}$, then the following inequality also holds
$
\left(
R_{1},w_{1}
\right)
\leq_{\mathbf{T}_{\Sigma}(X)}
\left(P,s
\right).
$

Now, consider the initial subsequence $(\mathfrak{R}_{k})_{k\in 2}$. This is a sequence of length $1$ instantiating that $(\mathfrak{Q},t)\prec_{\mathbf{Pth}_{\boldsymbol{\mathcal{A}}}} (\mathfrak{R}_{1},w_{1})$. Now, since $\mathfrak{R}_{1}$ is a $(1,0)$-identity  path, it follows, by the base step, that $\mathfrak{Q}$ is also a $(1,0)$-identity  path. 

Moreover, if the paths $\mathfrak{R}_{1}$ and $\mathfrak{Q}$ have the form
\begin{align*}
\mathfrak{Q}&=
\mathrm{ip}^{(1,0)\sharp}_{t}\left(
Q
\right),
&
\mathfrak{R}_{1}&=
\mathrm{ip}^{(1,0)\sharp}_{w_{1}}\left(
R_{1}
\right),
\end{align*}
for some terms $Q\in\mathrm{T}_{\Sigma}(X)_{t}$ and $R_{1}\in\mathrm{T}_{\Sigma}(X)_{w_{1}}$, then the following inequality also holds
$
\left(Q,t
\right)
\leq_{\mathbf{T}_{\Sigma}(X)}
\left(R_{1},w_{1}
\right).
$

All in all, we can affirm that $\mathfrak{Q}$ is a $(1,0)$-identity  path and, taking into account Remark~\ref{ROrd}, since $(\coprod\mathrm{T}_{\Sigma}(X),\leq_{\mathbf{T}_{\Sigma}(X)})$ is a partially ordered set, we have, by transitivity, that the following inequality holds
$
\left(Q,t
\right)
\leq_{\mathbf{T}_{\Sigma}(X)}
\left(P,s
\right).
$

This finishes the proof.
\end{proof}

From the last proposition we infer that every finite strictly decreasing sequence instantiating an inequality where the greatest element is a $(1,0)$-identity  path  is entirely composed of $(1,0)$-identity  paths. 

\begin{corollary}\label{CLower}
Let $m$ be a natural number in $\mathbb{N}-\{0\}$, $\mathbf{w}$ a word in $S^{\star}$ such that $\bb{\mathbf{w}}=m+1$, and $(\mathfrak{R}_{k})_{k\in \bb{\mathbf{w}}}$ a family of  paths in $\mathrm{Pth}_{\boldsymbol{\mathcal{A}},\mathbf{w}}$ such that
\begin{enumerate}
\item for every $k\in m$, $(\mathfrak{R}_{k}, w_{k})\prec_{\mathbf{Pth}_{\boldsymbol{\mathcal{A}}}} (\mathfrak{R}_{k+1}, w_{k+1})$; and 
\item $\mathfrak{R}_{m}$ is a $(1,0)$-identity path,
\end{enumerate}
then, for every $k\in m$, $\mathfrak{R}_{k}$ is a $(1,0)$-identity  path. Moreover, for every $k\in m+1$, if $\mathfrak{R}_{k}$ has the form
\begin{align*}
\mathfrak{R}_{k}&=
\mathrm{ip}^{(1,0)\sharp}_{w_{k}}\left(
R_{k}
\right),
\end{align*}
for some term $R_{w}\in\mathrm{T}_{\Sigma}(X)_{w_{k}}$ then, for every $k\in m$, the following strict inequality  holds
$$
\left(
R_{k},w_{k}\right)
<_{\mathbf{T}_{\Sigma}(X)}
\left(R_{k+1},w_{k+1}
\right).
$$
\end{corollary}

We next prove some technical results concerning the nature of pairs in the preorder $\leq_{\mathbf{Pth}_{\boldsymbol{\mathcal{A}}}}$. In this respect, we will show that, for every pair of sorts $s$ and $t$ in $S$ and every pair of paths $\mathfrak{P}\in\mathrm{Pth}_{\boldsymbol{\mathcal{A}},s}$, $\mathfrak{Q}\in\mathrm{Pth}_{\boldsymbol{\mathcal{A}},t}$, if $(\mathfrak{Q},t)<_{\mathbf{Pth}_{\boldsymbol{\mathcal{A}}}}(\mathfrak{P},s)$, then $\mathfrak{P}$ is strictly more complex than $\mathfrak{Q}$ with respect to either length  or  height of the associated translation.

\begin{lemma}\label{LOrdI} Let $m$ be a natural number in $\mathbb{N}-\{0\}$, $\mathbf{w}$ a word in $S^{\star}$ such that $\bb{\mathbf{w}}=m+1$, and $(\mathfrak{R}_{k})_{k\in \bb{\mathbf{w}}}$ a family of paths in $\mathrm{Pth}_{\boldsymbol{\mathcal{A}},\mathbf{w}}$ such that, for every $k\in m$, $(\mathfrak{R}_{k}, w_{k})\prec_{\mathbf{Pth}_{\boldsymbol{\mathcal{A}}}} (\mathfrak{R}_{k+1}, w_{k+1})$. Then, for every $k\in [1,m]$, we have that
\begin{enumerate}
\item $\mathfrak{R}_{k}$ is a $(1,0)$-identity path on a non-minimal term; or
\item $\mathfrak{R}_{k}$ is a path of length strictly greater than one containing echelons;  or
\item $\mathfrak{R}_{k}$ is an echelonless path.
\end{enumerate}
\end{lemma}
\begin{proof}
Let us recall (see Remark~\ref{RIrrefl}) that $\prec_{\mathbf{Pth}_{\boldsymbol{\mathcal{A}}}}$ and $\Delta_{\coprod\mathrm{Pth}_{\boldsymbol{\mathcal{A}}}}$ are disjoint sets. Moreover, by hypothesis, for every $k\in m$, $(\mathfrak{R}_{k}, w_{k})\prec_{\mathbf{Pth}_{\boldsymbol{\mathcal{A}}}} (\mathfrak{R}_{k+1}, w_{k+1})$. From this it follows that, for every $k\in [1,m]$, the element $(\mathfrak{R}_{k}, w_{k})$ is not a minimal element of $(\coprod \mathrm{Pth}_{\boldsymbol{\mathcal{A}}}, \leq_{\mathbf{Pth}_{\boldsymbol{\mathcal{A}}}})$.
\end{proof}

In the following lemma we show that, for every pair of sorts $s$ and $t$ in $S$ and every pair of paths 
$\mathfrak{P}\in\mathrm{Pth}_{\boldsymbol{\mathcal{A}},s}$ and $\mathfrak{Q}\in\mathrm{Pth}_{\boldsymbol{\mathcal{A}},t}$, if $(\mathfrak{Q},t)<_{\mathbf{Pth}_{\boldsymbol{\mathcal{A}}}}(\mathfrak{P},s)$ and there exists an strictly decreasing sequence instantiating this inequality  composed of paths of length strictly greater than one and containing each of them at least one echelon, then $\bb{\mathfrak{Q}}<\bb{\mathfrak{P}}$.

\begin{lemma}\label{LOrdII} Let $m$ be a natural number in $\mathbb{N}-\{0\}$, $\mathbf{w}$ a word in $S^{\star}$ such that $\bb{\mathbf{w}}=m+1$, and $(\mathfrak{R}_{k})_{k\in \bb{\mathbf{w}}}$ a family of paths in $\mathrm{Pth}_{\boldsymbol{\mathcal{A}},\mathbf{w}}$ such that 
\begin{enumerate}
\item for every $k\in m$, $(\mathfrak{R}_{k}, w_{k})\prec_{\mathbf{Pth}_{\boldsymbol{\mathcal{A}}}}(\mathfrak{R}_{k+1}, w_{k+1})$; and
\item for every $k\in [1,m]$, $\mathfrak{R}_{m}$ is a path of length strictly greater than one containing at least one echelon.
\end{enumerate}
Then $\bb{\mathfrak{R}_{0}}< \bb{\mathfrak{R}_{m}}$.
\end{lemma}
\begin{proof}
We prove it by induction on $m\in\mathbb{N}-\{0\}$.

\textsf{Base step of the induction.}

For $m=1$, we have that $(\mathfrak{R}_{0},w_{0})\prec_{\mathbf{Pth}_{\boldsymbol{\mathcal{A}}}}(\mathfrak{R}_{1},w_{1})$, where $\mathfrak{R}_{1}$ is a  path of length strictly greater than one containing at least one echelon.

Let $i\in\bb{\mathfrak{R}_{1}}$ be the first index for which $\mathfrak{R}^{i,i}_{1}$ is an echelon. Then, according to Definition~\ref{DOrd}, we have that $w_{0}=w_{1}$ and, taking into account the different possibilities for $i$, we have that $\mathfrak{R}_{0}$ has to be equal to one of the following subpaths of $\mathfrak{R}_{1}$
\begin{align*}
\mathfrak{R}^{0,0}_{1},&&
\mathfrak{R}^{1,\bb{\mathfrak{R}_{1}}-1}_{1},&&
\mathfrak{R}^{0,i-1}_{1}, \qquad \mbox{or}&&
\mathfrak{R}^{i,\bb{\mathfrak{R}_{1}}-1}_{1}.
\end{align*}
Note that, in every case, we can guarantee that  $\bb{\mathfrak{R}_{0}}<\bb{\mathfrak{R}_{1}}$.

\textsf{Inductive step of the induction.}

Let us suppose that the statement holds for sequences of length up to $m\in\mathbb{N}-\{0\}$. That is, if $\mathbf{w}$ is a word in $S^{\star}$ of length $\bb{\mathbf{w}}=m+1$ and $(\mathfrak{R}_{k})_{k\in\bb{\mathbf{w}}}$ is a family of  paths in $\mathrm{Pth}_{\boldsymbol{\mathcal{A}},\mathbf{w}}$ satisfying that 
\begin{enumerate}
\item for every $k\in m$, $(\mathfrak{R}_{k}, w_{k})\prec_{\mathbf{Pth}_{\boldsymbol{\mathcal{A}}}} (\mathfrak{R}_{k+1}, w_{k+1})$; and
\item for every $k\in [1,m]$, $\mathfrak{R}_{k}$ is a  path of length strictly greater than one containing at least one echelon,
\end{enumerate}
then
$\bb{\mathfrak{R}_{0}}<\bb{\mathfrak{R}_{m}}$.

Now, let $\mathbf{w}$ be a word in $S^{\star}$ of length $\bb{\mathbf{w}}=m+2$ and let $(\mathfrak{R}_{k})_{k\in\bb{\mathbf{w}}}$ be a family of  paths in $\mathrm{Pth}_{\boldsymbol{\mathcal{A}},\mathbf{w}}$ satisfying that 
\begin{enumerate}
\item for every $k\in m+1$, $(\mathfrak{R}_{k}, w_{k})\prec_{\mathbf{Pth}_{\boldsymbol{\mathcal{A}}}} (\mathfrak{R}_{k+1}, w_{k+1})$; and
\item for every $k\in [1,m+1]$, $\mathfrak{R}_{k}$ is a path of length strictly greater than one containing at least one echelon.
\end{enumerate}

Consider the final subsequence $(\mathfrak{R}_{k})_{k\in [1,m+1]}$. This is a sequence of length $m+1$ satisfying the conditions of the inductive step. Hence, $\bb{\mathfrak{R}_{1}}<\bb{\mathfrak{R}_{m+1}}$.

Now, consider the initial subsequence $(\mathfrak{R}_{k})_{k\in 2}$. This is a sequence of length $2$ instantiating that $(\mathfrak{R}_{0},w_{0})\prec_{\mathbf{Pth}_{\boldsymbol{\mathcal{A}}}}(\mathfrak{R}^{(2)}_{1},w_{1})$ and, by assumption, $\mathfrak{R}_{1}$ is a  path of length strictly greater than one containing at least one  echelon. By the base step, we have that $\bb{\mathfrak{R}_{0}}<\bb{\mathfrak{R}_{1}}$.

All in all, we can affirm that $\bb{\mathfrak{R}_{0}}<\bb{\mathfrak{R}_{m+1}}$.

This finishes the proof.
\end{proof}

In the following lemma we show that, for every pair of sorts $s$ and $t$ in $S$ and every pair of  paths $\mathfrak{P}\in\mathrm{Pth}_{\boldsymbol{\mathcal{A}},s}$ and $\mathfrak{Q}\in\mathrm{Pth}_{\boldsymbol{\mathcal{A}},t}$, if $(\mathfrak{Q},t)<_{\mathbf{Pth}_{\boldsymbol{\mathcal{A}}}}(\mathfrak{P},s)$ and there exists an strictly decreasing sequence instantiating this inequality composed of echelonless paths, then $\bb{\mathfrak{Q}}\leq \bb{\mathfrak{P}}$ and if $\mathfrak{Q}$ is a path of length at least one, then the maximum of the heights of the translations occurring in $\mathfrak{Q}$ is strictly smaller than the maximum of the height of the translations occurring in $\mathfrak{P}$.

\begin{lemma}\label{LOrdIII} 
Let $m$ be a natural number in $\mathbb{N}-\{0\}$, $\mathbf{w}$ a word in $S^{\star}$ such that $\bb{\mathbf{w}}=m+1$, and  $(\mathfrak{R}_{k})_{k\in \bb{\mathbf{w}}}$ a family of  paths in $\mathrm{Pth}_{\boldsymbol{\mathcal{A}},\mathbf{w}}$ such that
\begin{enumerate}
\item for every $k\in m$, $(\mathfrak{R}_{k}, w_{k})\prec_{\mathbf{Pth}_{\boldsymbol{\mathcal{A}}}} (\mathfrak{R}_{k+1}, w_{k+1})$; and
\item for every $k\in [1,m]$, $\mathfrak{R}_{k}$ is an echelonless path.
\end{enumerate}
Then $\bb{\mathfrak{R}_{0}}\leq \bb{\mathfrak{R}_{m}}$. Moreover, if $\bb{\mathfrak{R}_{0}}\geq 1$, and the  paths $\mathfrak{R}_{0}$ and $\mathfrak{R}_{m}$ have the form
\begin{align*}
\mathfrak{R}_{0}&=\left(
(R_{0,j})_{j\in\bb{\mathbf{c}}+1},
(\mathfrak{r}_{0,j})_{j\in\bb{\mathbf{c}}},
(T_{0,j})_{j\in\bb{\mathbf{c}}}
\right)
\\
\mathfrak{R}_{m}&=\left(
(R_{m,j})_{j\in\bb{\mathbf{d}}+1},
(\mathfrak{r}_{m,j})_{j\in\bb{\mathbf{d}}},
(T_{m,j})_{j\in\bb{\mathbf{d}}}
\right),
\end{align*}
for some non-empty words $\mathbf{c},\mathbf{d}\in S^{\star}-\{\lambda\}$, then 
\[
\max\{
\bb{T_{0,j}}\mid j\in \bb{\mathbf{c}}
\}
<
\max\{
\bb{T_{m,j}}\mid j\in \bb{\mathbf{d}}
\}.
\]
\end{lemma}
\begin{proof}
We prove it by induction on $m\in\mathbb{N}-\{0\}$.
 
\textsf{Base step of the induction}. 

For $m=1$ we have that $(\mathfrak{R}_{0},w_{0})\prec_{\mathbf{Pth}_{\boldsymbol{\mathcal{A}}}}(\mathfrak{R}_{1},w_{1})$, where $\mathfrak{R}_{1}$ is an echelonless  $\mathbf{d}$-path in $\mathrm{Pth}_{\boldsymbol{\mathcal{A}},w_{1}}$, for some non-empty word $\mathbf{d}\in S^{\star}-\{\lambda\}$, of the form
$$
\mathfrak{R}_{1}=
\left(
(R_{1,i})_{i\in\bb{\mathbf{d}}+1},
(\mathfrak{p}_{1,i})_{i\in\bb{\mathbf{d}}},
(T_{1,i})_{i\in\bb{\mathbf{d}}}
\right).
$$
such that, for a unique word $\mathbf{w_{1}}\in S^{\star}-\{\lambda\}$, a unique operation symbol $\sigma\in\Sigma_{\mathbf{w_{1}},w_{1}}$, 
$(T_{1,i})_{i\in\bb{\mathbf{d}}}$ is a family of translations of type $\sigma$.

Let $((\mathbf{d}_{j})_{j\in\bb{\mathbf{w_{1}}}}, (\mathfrak{R}_{1,j})_{j\in\bb{\mathbf{w_{1}}}})$ be the pair in $(S^{\star})^{\bb{\mathbf{w_{1}}}}\times \mathrm{Pth}_{\boldsymbol{\mathcal{A}},\mathbf{w_{1}}}$ we can extract from $\mathfrak{R}_{1}$ after applying the extraction path algorithm from Lemma~\ref{LPthExtract}. Since $(\mathfrak{R}_{0},w_{0})\prec_{\mathbf{Pth}_{\boldsymbol{\mathcal{A}}}}(\mathfrak{R}_{1},w_{1})$, according to Definition~\ref{DOrd}, there exists some $j\in\bb{\mathbf{w_{1}}}$ for which $w_{0}=\mathbf{w}_{\mathbf{1},j}$ and $\mathfrak{R}_{0}=\mathfrak{R}_{1,j}$. Let us recall from the proof of Lemma~\ref{LPthExtract}, that $\sum_{j\in\bb{\mathbf{w_{1}}}}\bb{\mathbf{d}_{j}}=\bb{\mathbf{d}}$. Consequently, 
$$\bb{\mathfrak{R}_{0}}=\bb{\mathbf{d}_{j}}\leq\bb{\mathbf{d}}=\bb{\mathfrak{R}_{1}}.$$ 

Moreover, if $\bb{\mathfrak{R}_{0}}\geq 1$ and $\mathfrak{R}_{0}$ has the form
$$\mathfrak{R}_{0}=\left(
(R_{0,j})_{j\in\bb{\mathbf{d}_{j}}+1},
(\mathfrak{r}_{0,j})_{j\in\bb{\mathbf{d}_{j}}},
(T_{0,j})_{j\in\bb{\mathbf{d}_{j}}}
\right),$$
we have, according to the proof of Lemma~\ref{LPthExtract}, that $(T_{0,j})_{j\in\bb{\mathbf{d}_{j}}}$ is a family of derived  translations of a subsequence of $(T_{1,i})_{i\in\bb{\mathbf{d}}}$, hence all the  translations occurring in $\mathfrak{R}_{0}$ have their heights strictly bounded by the height of a  translation occurring in $\mathfrak{R}_{1}$, hence
\[
\max\{
\bb{T_{0,j}}\mid j\in\bb{\mathbf{d}_{j}}
\}
<
\max\{
\bb{T_{1,i}}\mid i\in \bb{\mathbf{d}}
\}.
\]

\textsf{Inductive step of the induction}. 

Let us suppose that the statement holds for sequences of length up to $m\in\mathbb{N}-\{0\}$. That is, if $\mathbf{w}$ is a word in $S^{\star}$ of length $\bb{\mathbf{w}}=m+1$ and  $(\mathfrak{R}_{k})_{k\in \bb{\mathbf{w}}}$ a family of  paths in $\mathrm{Pth}_{\boldsymbol{\mathcal{A}},\mathbf{w}}$ satisfying
\begin{enumerate}
\item for every $k\in m$, $(\mathfrak{R}_{k}, w_{k})\prec_{\mathbf{Pth}_{\boldsymbol{\mathcal{A}}}} (\mathfrak{R}_{k+1}, w_{k+1})$; and
\item for every $k\in [1,m]$, $\mathfrak{R}_{k}$ is an echelonless path,
\end{enumerate}
then $\bb{\mathfrak{R}_{0}}\leq \bb{\mathfrak{R}_{m}}$.

Moreover, if $\bb{\mathfrak{R}_{0}}\geq 1$, and the  paths $\mathfrak{R}_{0}$ and $\mathfrak{R}_{m}$ have the form
\begin{align*}
\mathfrak{R}_{0}&=\left(
(R_{0,j})_{j\in\bb{\mathbf{c}}+1},
(\mathfrak{r}_{0,j})_{j\in\bb{\mathbf{c}}},
(T_{0,j})_{j\in\bb{\mathbf{c}}}
\right)
\\
\mathfrak{R}_{m}&=\left(
(R_{m,j})_{j\in\bb{\mathbf{d}}+1},
(\mathfrak{r}_{m,j})_{j\in\bb{\mathbf{d}}},
(T_{m,j})_{j\in\bb{\mathbf{d}}}
\right),
\end{align*}
for some non-empty words $\mathbf{c},\mathbf{d}\in S^{\star}-\{\lambda\}$, then 
\[
\max\{
\bb{T_{0,j}}\mid j\in\bb{\mathbf{c}}
\}
<
\max\{
\bb{T_{m,j}}\mid j\in \bb{\mathbf{d}}
\}.
\]

Now, let $\mathbf{w}$ be a word in $S^{\star}$ of length $\bb{\mathbf{w}}=m+2$ and $(\mathfrak{R}_{k})_{k\in \bb{\mathbf{w}}}$ a family of  paths in $\mathrm{Pth}_{\boldsymbol{\mathcal{A}},\mathbf{w}}$ satisfying
\begin{enumerate}
\item for every $k\in m+1$, $(\mathfrak{R}_{k}, w_{k})\prec_{\mathbf{Pth}_{\boldsymbol{\mathcal{A}}}} (\mathfrak{R}_{k+1}, w_{k+1})$; and
\item for every $k\in [1,m+1]$, $\mathfrak{R}_{k}$ is an echelonless  path.
\end{enumerate}

Consider the final subsequence $(\mathfrak{R}_{k})_{k\in[1,m+1]}$. This is a sequence of length $m+1$ satisfying the conditions of the inductive step. Hence, 
$\bb{\mathfrak{R}_{1}}\leq\bb{\mathfrak{R}_{m+1}}$. 

Moreover, by assumption, $\mathfrak{R}_{1}$ is an echelonless  path. Thus, if $\mathfrak{R}_{1}$ and $\mathfrak{R}_{m+1}$ have the form
\begin{align*}
\mathfrak{R}_{1}&=\left(
(R_{1,j})_{j\in\bb{\mathbf{d}}+1},
(\mathfrak{r}_{1,j})_{j\in\bb{\mathbf{d}}},
(T_{1,j})_{j\in\bb{\mathbf{d}}}
\right)
\\
\mathfrak{R}_{m+1}&=\left(
(R_{m+1,j})_{j\in\bb{\mathbf{e}}+1},
(\mathfrak{r}_{m+1,j})_{j\in\bb{\mathbf{e}}},
(T_{m+1,j})_{j\in\bb{\mathbf{e}}}
\right),
\end{align*}
for some non-empty words $\mathbf{d},\mathbf{e}\in S^{\star}-\{\lambda\}$, then 
\[
\max\{
\bb{T_{1,j}}\mid j\in\bb{\mathbf{d}}
\}
<
\max\{
\bb{T_{m+1,j}}\mid j\in \bb{\mathbf{e}}
\}.
\]

Now, consider the initial subsequence $(\mathfrak{R}_{k})_{k\in 2}$. This is a sequence of length $2$ instantiating that $(\mathfrak{R}_{0},w_{0})\prec_{\mathbf{Pth}_{\boldsymbol{\mathcal{A}}}}(\mathfrak{R}_{1},w_{1})$ and, by assumption, $\mathfrak{R}_{1}$ is an echelonless  path. By the base step, we have that $\bb{\mathfrak{R}_{0}}\leq \bb{\mathfrak{R}_{1}}$.

All in all, we can affirm that $\bb{\mathfrak{R}_{0}}\leq\bb{\mathfrak{R}_{m+1}}$.

Moreover, if $\bb{\mathfrak{R}_{0}}\geq 1$ and the  path $\mathfrak{R}_{0}$ has the form
\begin{align*}
\mathfrak{R}_{0}&=\left(
(R_{0,j})_{j\in\bb{\mathbf{c}}+1},
(\mathfrak{r}_{0,j})_{j\in\bb{\mathbf{c}}},
(T_{0,j})_{j\in\bb{\mathbf{c}}}
\right),
\end{align*}
for some non-empty word $\mathbf{c}\in S^{\star}-\{\lambda\}$, then 
\[
\max\{
\bb{T_{0,j}}\mid j\in\bb{\mathbf{c}}
\}
<
\max\{
\bb{T_{1,j}}\mid j\in\bb{\mathbf{d}}
\}.
\]

All in all, we can affirm in this case that 
\[
\max\{
\bb{T_{0,j}}\mid j\in\bb{\mathbf{c}}
\}
<
\max\{
\bb{T_{m+1,j}}\mid j\in \bb{\mathbf{e}}
\}.
\]

This finishes the proof.
\end{proof}

In the following corollary we show that, for every pair of sorts $s$ and $t$ in $S$ and every pair of paths $\mathfrak{P}\in\mathrm{Pth}_{\boldsymbol{\mathcal{A}},s}$ and $\mathfrak{Q}\in\mathrm{Pth}_{\boldsymbol{\mathcal{A}},t}$, if $(\mathfrak{Q},t)<_{\mathbf{Pth}_{\boldsymbol{\mathcal{A}}}}(\mathfrak{P},s)$ and $\mathfrak{P}$ is a path of length strictly greater than one containing at least one  echelon then $\bb{\mathfrak{Q}}<\bb{\mathfrak{P}}$.

\begin{corollary}\label{COrdI}Let $m$ be a natural number in $\mathbb{N}-\{0\}$, $\mathbf{w}$ a word in $S^{\star}$ such that $\bb{\mathbf{w}}=m+1$, and  $(\mathfrak{R}_{k})_{k\in \bb{\mathbf{w}}}$ a family of  paths in $\mathrm{Pth}_{\boldsymbol{\mathcal{A}},\mathbf{w}}$ such that
\begin{enumerate}
\item  for every $k\in m$, $(\mathfrak{R}_{k}, w_{k})\prec_{\mathbf{Pth}_{\boldsymbol{\mathcal{A}}}} (\mathfrak{R}_{k+1}, w_{k+1})$; and
\item $\mathfrak{R}_{m}$ is a  path of length strictly greater than one containing at least one  echelon.
\end{enumerate}
Then $\bb{\mathfrak{R}_{0}}<\bb{\mathfrak{R}_{m}}$.
\end{corollary}
\begin{proof}
We prove it by induction on $m\in\mathbb{N}-\{0\}$.

\textsf{Base step of the induction.}

For $m=1$, we have that $(\mathfrak{R}_{0},w_{0})\prec_{\mathbf{Pth}_{\boldsymbol{\mathcal{A}}}}(\mathfrak{R}_{1},w_{1})$, where $\mathfrak{R}_{1}$ is a path of length strictly greater than one containing at least one echelon. Then by  Lemma~\ref{LOrdII}, we have that 
$\bb{\mathfrak{R}_{0}}<\bb{\mathfrak{R}_{1}}$.

\textsf{Inductive step of the induction.}

Let us suppose that the statement holds for sequences of length up to $m\in\mathbb{N}-\{0\}$. That is, if $\mathbf{w}$ is a word in $S^{\star}$ of length $\bb{\mathbf{w}}=m+1$ and $(\mathfrak{R}_{k})_{k\in \bb{\mathbf{w}}}$ a family of paths in $\mathrm{Pth}_{\boldsymbol{\mathcal{A}},\mathbf{w}}$ such that
\begin{enumerate}
\item  for every $k\in m$, $(\mathfrak{R}_{k}, w_{k})\prec_{\mathbf{Pth}_{\boldsymbol{\mathcal{A}}}} (\mathfrak{R}_{k+1}, w_{k+1})$; and
\item $\mathfrak{R}_{m}$ is a path of length strictly greater than one containing at least one echelon,
\end{enumerate}
then $\bb{\mathfrak{R}_{0}}<\bb{\mathfrak{R}_{m}}$.

Now, let $\mathbf{w}$ be a word in $S^{\star}$ of length $\bb{\mathbf{w}}=m+2$ and $(\mathfrak{R}_{k})_{k\in \bb{\mathbf{w}}}$ a family of paths in $\mathrm{Pth}_{\boldsymbol{\mathcal{A}},\mathbf{w}}$ such that
\begin{enumerate}
\item  for every $k\in m+1$, $(\mathfrak{R}_{k}, w_{k})\prec_{\mathbf{Pth}_{\boldsymbol{\mathcal{A}}}} (\mathfrak{R}_{k+1}, w_{k+1})$; and
\item $\mathfrak{R}_{m+1}$ is a path of length strictly greater than one containing at least one echelon.
\end{enumerate}

Consider the final subsequence $(\mathfrak{R}_{k})_{k\in[1,m+1]}$. This is a sequence of length $m+1$ satisfying the conditions of the inductive step. Hence $\bb{\mathfrak{R}_{1}}<\bb{\mathfrak{R}_{m+1}}$.

Now, consider the initial subsequence $(\mathfrak{R}_{k})_{k\in 2}$. This is a sequence of length $2$ instantiating that $(\mathfrak{R}_{0},w_{0})\prec_{\mathbf{Pth}_{\boldsymbol{\mathcal{A}}}} (\mathfrak{R}_{1},w_{1})$. Since $\prec_{\mathbf{Pth}_{\boldsymbol{\mathcal{A}}}}$ is an irreflexive relation by Remark~\ref{RIrrefl}, we conclude that $(\mathfrak{R}_{1},w_{1})$ cannot be 
a minimal element of $\leq_{\mathbf{Pth}_{\boldsymbol{\mathcal{A}}}}$. Thus, taking into account Lemma~\ref{LOrdI}, the different possibilities for $\mathfrak{R}_{1}$ are
\begin{enumerate}
\item[(1)] $\mathfrak{R}_{1}$ is a $(1,0)$-identity  path on a non-minimal term. Following Proposition~\ref{PLower}, we conclude that $\mathfrak{R}_{0}$ is also a $(1,0)$-identity path. Hence, $\bb{\mathfrak{R}_{0}}=\bb{\mathfrak{R}_{1}}$; or
\item[(2)] $\mathfrak{R}_{1}$ is a path of length strictly greater than one containing at least one echelon. By Lemma~\ref{LOrdII} we conclude that $\bb{\mathfrak{R}_{0}}<\bb{\mathfrak{R}_{1}}$; or
\item[(3)] $\mathfrak{R}_{1}$ is an echelonless  path. By Lemma~\ref{LOrdIII} we conclude that $\bb{\mathfrak{R}_{0}}\leq\bb{\mathfrak{R}_{1}}$.
\end{enumerate}
All in all, we conclude that 
$\bb{\mathfrak{R}_{0}}<\bb{\mathfrak{R}_{m+1}}$.

This finishes the proof.
\end{proof}

In the following corollary we show that, for every pair of sorts $s$ and $t$ in $S$ and every pair of paths $\mathfrak{P}\in\mathrm{Pth}_{\boldsymbol{\mathcal{A}},s}$ and $\mathfrak{Q}\in\mathrm{Pth}_{\boldsymbol{\mathcal{A}},t}$, if $(\mathfrak{Q},t)<_{\mathbf{Pth}_{\boldsymbol{\mathcal{A}}}}(\mathfrak{P},s)$ and $\mathfrak{P}$ is either (1) a  path of length strictly greater than one containing at least one  echelon; or (2) an echelonless path, then $\bb{\mathfrak{Q}}\leq \bb{\mathfrak{P}}$.

\begin{corollary}\label{COrdII}Let $m$ be a natural number in $\mathbb{N}-\{0\}$, $\mathbf{w}$ a word in $S^{\star}$ such that $\bb{\mathbf{w}}=m+1$, and  $(\mathfrak{R}_{k})_{k\in \bb{\mathbf{w}}}$ a family of paths in $\mathrm{Pth}_{\boldsymbol{\mathcal{A}},\mathbf{w}}$ such that
\begin{enumerate}
\item for every $k\in m$, $(\mathfrak{R}_{k}, w_{k})\prec_{\mathbf{Pth}_{\boldsymbol{\mathcal{A}}}} (\mathfrak{R}_{k+1}, w_{k+1})$;
\end{enumerate}
and either one of the following conditions holds
\begin{enumerate}
\item[(2.1)] $\mathfrak{R}_{m}$ is a path of length strictly greater than one containing at least one  echelon; or 
\item[(2.2)] $\mathfrak{R}_{m}$ an echelonless path.
\end{enumerate}
Then $\bb{\mathfrak{R}_{0}}\leq \bb{\mathfrak{R}_{m}}$.
\end{corollary}
\begin{proof}
We prove it by induction on $m\in\mathbb{N}-\{0\}$.

\textsf{Base step of the induction.}

For $m=1$, we have that  $(\mathfrak{R}_{0},w_{0})\prec_{\mathbf{Pth}_{\boldsymbol{\mathcal{A}}}}(\mathfrak{R}_{1},w_{1})$, where either (2.1) $\mathfrak{R}_{1}$ is a path of length strictly greater than one containing at least one echelon; or  $\mathfrak{R}_{1}$ an echelonless  path. Then by either Lemma~\ref{LOrdII} or~\ref{LOrdIII}, respectively, we have that $\bb{\mathfrak{R}_{0}}\leq\bb{\mathfrak{R}_{1}}$.

\textsf{Inductive step of the induction.}

Let us suppose that the statement holds for sequences of length up to $m\in\mathbb{N}-\{0\}$. That is, if 
$\mathbf{w}$ is a word in $S^{\star}$ of length $\bb{\mathbf{w}}=m+1$ and $(\mathfrak{R}_{k})_{k\in\bb{\mathbf{w}}}$ a family of paths in $\mathrm{Pth}_{\boldsymbol{\mathcal{A}},\mathbf{w}}$ such that
\begin{enumerate}
\item for every $k\in m$, $(\mathfrak{R}_{k}, w_{k})\prec_{\mathbf{Pth}_{\boldsymbol{\mathcal{A}}}} (\mathfrak{R}_{k+1}, w_{k+1})$;
\end{enumerate}
and  either one of the following conditions holds
\begin{enumerate}
\item[(2.1)] $\mathfrak{R}_{m}$ is a path of length strictly greater than one containing at least one echelon; or 
\item[(2.2)] $\mathfrak{R}_{m}$ is an echelonless path.
\end{enumerate}
Then $\bb{\mathfrak{R}_{0}}\leq \bb{\mathfrak{R}_{m}}$.

Now let $\mathbf{w}$ be a word in $S^{\star}$ of length $\bb{\mathbf{w}}=m+2$ and let $(\mathfrak{R}_{k})_{k\in\bb{\mathbf{w}}}$ be a family of  paths in $\mathrm{Pth}_{\boldsymbol{\mathcal{A}},\mathbf{w}}$ such that
\begin{enumerate}
\item for every $k\in m+1$, $(\mathfrak{R}_{k}, w_{k})\prec_{\mathbf{Pth}_{\boldsymbol{\mathcal{A}}}} (\mathfrak{R}_{k+1}, w_{k+1})$;
\end{enumerate}
and  either one of the following conditions holds
\begin{enumerate}
\item[(2.1)] $\mathfrak{R}_{m+1}$ is a  path of length strictly greater than one containing at least one  echelon; or 
\item[(2.2)] $\mathfrak{R}_{m+1}$ is an echelonless path.
\end{enumerate}

Consider the final subsequence $(\mathfrak{R}_{k})_{k\in[1,m+1]}$. This is a sequence of length $m+1$ satisfying the conditions of the inductive step. Hence $\bb{\mathfrak{R}_{1}}\leq\bb{\mathfrak{R}_{m+1}}$.

Now, consider the initial subsequence $(\mathfrak{R}_{k})_{k\in 2}$. This is a sequence of length $2$ instantiating that $(\mathfrak{R}_{0},w_{0})\prec_{\mathbf{Pth}_{\boldsymbol{\mathcal{A}}}} (\mathfrak{R}_{1},w_{1})$. Since $\prec_{\mathbf{Pth}_{\boldsymbol{\mathcal{A}}}}$ is an irreflexive relation by Remark~\ref{RIrrefl}, we conclude that $(\mathfrak{R}_{1},w_{1})$ cannot be a minimal element of 
$\leq_{\mathbf{Pth}_{\boldsymbol{\mathcal{A}}}}$. Thus, taking into account Lemma~\ref{LOrdI}, the different possibilities for $\mathfrak{R}_{1}$ are
\begin{enumerate}
\item[(1)] $\mathfrak{R}_{1}$ is a $(1,0)$-identity path on a non-minimal term. Following Proposition~\ref{PLower}, we conclude that $\mathfrak{R}_{0}$ is also a $(1,0)$-identity path. Hence, $\bb{\mathfrak{R}_{0}}=\bb{\mathfrak{R}_{1}}$; or
\item[(2)] $\mathfrak{R}_{1}$ is a path of length strictly greater than one containing at least one  echelon. By Lemma~\ref{LOrdII} we conclude that 
$\bb{\mathfrak{R}_{0}}<\bb{\mathfrak{R}_{1}}$; or
\item[(3)] $\mathfrak{R}_{1}$ is an echelonless path. Then by Lemma~\ref{LOrdIII} we have that $\bb{\mathfrak{R}_{0}}\leq\bb{\mathfrak{R}_{1}}$.
\end{enumerate}
All in all, we conclude that $\bb{\mathfrak{R}_{0}}\leq\bb{\mathfrak{R}_{m+1}}$.

This finishes the proof.
\end{proof}

The above technical lemmas will be used to prove the following proposition.

\begin{restatable}{proposition}{POrdArt}
\label{POrdArt}
The preorder $\leq_{\mathbf{Pth}_{\boldsymbol{\mathcal{A}}}}$ on $\coprod\mathrm{Pth}_{\boldsymbol{\mathcal{A}}}$ is antisymmetric and, then, in the ordered set $(\coprod\mathrm{Pth}_{\boldsymbol{\mathcal{A}}}, \leq_{\mathbf{Pth}_{\boldsymbol{\mathcal{A}}}})$ there is not any strictly decreasing $\omega_{0}$-chain, i.e., $(\coprod\mathrm{Pth}_{\boldsymbol{\mathcal{A}}}, \leq_{\mathbf{Pth}_{\boldsymbol{\mathcal{A}}}})$ is an Artinian ordered set.
\end{restatable}
\begin{proof}
We first prove that $\leq_{\mathbf{Pth}_{\boldsymbol{\mathcal{A}}}}$ is antisymmetric. 

Let $s$ and $t$ be sorts in $S$, $\mathfrak{P}\in\mathrm{Pth}_{\boldsymbol{\mathcal{A}},s}$, and $\mathfrak{Q}\in\mathrm{Pth}_{\boldsymbol{\mathcal{A}},t}$. Let us suppose that $(\mathfrak{Q},t)\leq_{\mathbf{Pth}_{\boldsymbol{\mathcal{A}}}} (\mathfrak{P},s)$ and $(\mathfrak{P},s)\leq_{\mathbf{Pth}_{\boldsymbol{\mathcal{A}}}} (\mathfrak{Q},t)$. We want to show that $s=t$ and $\mathfrak{Q}=\mathfrak{P}$.

This statement is trivial if $\mathfrak{P}$ or $\mathfrak{Q}$ is a minimal element of $(\coprod\mathrm{Pth}_{\boldsymbol{\mathcal{A}}},\leq_{\mathbf{Pth}_{\boldsymbol{\mathcal{A}}}})$. So let us assume that neither $\mathfrak{P}$ nor $\mathfrak{Q}$ are minimal in $(\coprod\mathrm{Pth}_{\boldsymbol{\mathcal{A}}},\leq_{\mathbf{Pth}_{\boldsymbol{\mathcal{A}}}})$. From Lemma~\ref{LOrdI}, for $\mathfrak{P}$ we have that either (A) it is a $(1,0)$-identity path on a non-minimal term, or (B) it is a path of length strictly greater than one containing at least one echelon or (C) it is an echelonless path. In the same way, for $\mathfrak{Q}$ we have that either (D) it is a $(1,0)$-identity path on a non-minimal term, or (E) it is a  path of length strictly greater than one containing at least one echelon or (F) it is an echelonless path. However, as we will show immediately below only cases $(A,D)$, $(B,E)$ and $(C,F)$ are feasible.

We first consider the case (A) where $\mathfrak{P}$ is a $(1,0)$-identity path. Since we are assuming that $(\mathfrak{Q},t)\leq_{\mathbf{Pth}_{\boldsymbol{\mathcal{A}}}} (\mathfrak{P},s)$ we conclude, by Proposition~\ref{PLower}, that $\mathfrak{Q}$ is a $(1,0)$-identity path, i.e., (D). Hence, the $(1,0)$-identity paths $\mathfrak{P}$ and $\mathfrak{Q}$ have the form
\begin{align*}
\mathfrak{Q}&=\mathrm{ip}^{(1,0)\sharp}_{t}\left(
Q
\right),
&
\mathfrak{P}&=\mathrm{ip}^{(1,0)\sharp}_{s}\left(
P
\right),
\end{align*}
for a suitable pair of terms $P\in\mathrm{T}_{\Sigma}(X)_{s}$ and $Q\in\mathrm{T}_{\Sigma}(X)_{t}$. But, since we are assuming that $(\mathfrak{Q},t)\leq_{\mathbf{Pth}_{\boldsymbol{\mathcal{A}}}} (\mathfrak{P},s)$ and 
$(\mathfrak{P},s)\leq_{\mathbf{Pth}_{\boldsymbol{\mathcal{A}}}} (\mathfrak{Q},t)$, we conclude that 
$(Q,t)\leq_{\mathbf{T}_{\Sigma}(X)} (P,s)$ and $(P,t)\leq_{\mathbf{T}_{\Sigma}(X)} (Q,s)$. However, since 
$(\coprod\mathrm{T}_{\Sigma}(X), \leq_{\mathbf{T}_{\Sigma}(X)})$, by Remark~\ref{ROrd}, is an ordered set,   we conclude that $(Q,t)=(P,s)$. In particular $s=t$ and 
$$\mathfrak{Q}=\mathrm{ip}^{(1,0)\sharp}_{t}\left(
Q
\right)=\mathrm{ip}^{(1,0)\sharp}_{s}\left(
P
\right)=\mathfrak{P}.$$

The same result is obtained assuming that (D), i.e., that $\mathfrak{Q}$ is a $(1,0)$-identity path. Therefore cases (A) and (D) are totally determined.

We next show that cases (B,F) and (C,E) are infeasible.

If (B,F), then $\mathfrak{P}$ and $\mathfrak{Q}$ are necessarily different since one contains at least one  echelon and the other is echelonless. But $(\mathfrak{Q},t)<_{\mathbf{Pth}_{\boldsymbol{\mathcal{A}}}} (\mathfrak{P},s)$ and $\mathfrak{P}$ is a path of length strictly greater than one containing at least one echelon, hence, by Corollary~\ref{COrdI}, we have that $\bb{\mathfrak{Q}}<\bb{\mathfrak{P}}$. Moreover, since $(\mathfrak{P},s)<_{\mathbf{Pth}_{\boldsymbol{\mathcal{A}}}} (\mathfrak{Q},t)$ and $\mathfrak{Q}$ is an echelonless path, we have, by Corollary~\ref{COrdII}, that $\bb{\mathfrak{P}}\leq \bb{\mathfrak{Q}}$. This leads to the following contradiction
$$
\bb{\mathfrak{P}}\leq\bb{\mathfrak{Q}}<\bb{\mathfrak{P}}.
$$

The same conclusion can be drawn for case (C,E), by applying to it the same reasoning as in the case (B,F).

Therefore we are only left with these options: (B,E) and (C,F), where the paths $\mathfrak{P}$ and $\mathfrak{Q}$ are required to have the same nature.

For (B,E), let us suppose, towards a contradiction, that $\mathfrak{P}\neq\mathfrak{Q}$. Then, since $(\mathfrak{Q},t)<_{\mathbf{Pth}_{\boldsymbol{\mathcal{A}}}}(\mathfrak{P},s)$ and $\mathfrak{P}$ is a path of length strictly greater than one containing at least one  echelon, hence, by Corollary~\ref{COrdI}, we have that $\bb{\mathfrak{Q}}<\bb{\mathfrak{P}}$. Moreover, since $(\mathfrak{P},s)<_{\mathbf{Pth}_{\boldsymbol{\mathcal{A}}}}(\mathfrak{Q},t)$ and $\mathfrak{Q}$  is a path of length strictly greater than one containing at least one  echelon, hence, by Corollary~\ref{COrdI}, we have that $\bb{\mathfrak{P}}<\bb{\mathfrak{Q}}$. This leads to the following contradiction
$$
\bb{\mathfrak{P}}<\bb{\mathfrak{Q}}<\bb{\mathfrak{P}}.
$$
Therefore, we infer that $\mathfrak{P}=\mathfrak{Q}$.

For (C,F), let us suppose, towards a contradiction, that $\mathfrak{P}\neq \mathfrak{Q}$. Then, since $(\mathfrak{Q},t)<_{\mathbf{Pth}_{\boldsymbol{\mathcal{A}}}} (\mathfrak{P},s)$, we have, by Remark~\ref{ROrd}, that there exists a natural number $m\in\mathbb{N}-\{0\}$, a word $\mathbf{w}\in S^{\star}$ of length $\bb{\mathbf{w}}=m+1$, and a family of paths $(\mathfrak{R}_{k})_{k\in\bb{\mathbf{w}}}$ in $\mathrm{Pth}_{\boldsymbol{\mathcal{A}},\mathbf{w}}$, such that $w_{0}=t$, $\mathfrak{R}_{0}=\mathfrak{Q}$, $w_{m}=s$, $\mathfrak{R}_{m}=\mathfrak{P}$ and, for every $k\in m$, $(\mathfrak{R}_{k},w_{k})\prec_{\mathbf{Pth}_{\boldsymbol{\mathcal{A}}}} (\mathfrak{R}_{k+1},w_{k+1})$. Then, by Lemma~\ref{LOrdI}, for every $k\in [1,m]$, we have that 
\begin{enumerate}
\item $\mathfrak{R}_{k}$ is a $(1,0)$-identity path on a non-minimal term; or
\item $\mathfrak{R}_{k}$ is a path of length strictly greater than one containing at least one echelon; or 
\item $\mathfrak{R}_{k}$ is an echelonless path.
\end{enumerate}

We claim that, among all the above possibilities, for every $k\in [1,m]$, the path $\mathfrak{R}_{k}$ has to be an echelonless path. 

Let us suppose, towards a contradiction, that there exists an index $k\in [1,m]$ for which $\mathfrak{R}_{k}$ is a $(1,0)$-identity  path on a non-minimal term. Then, taking into account that $(\mathfrak{Q},t)\leq_{\mathbf{Pth}_{\boldsymbol{\mathcal{A}}}} (\mathfrak{R}_{k},w_{k})$, we obtain, 
according to Proposition~\ref{PLower}, that $\mathfrak{Q}$ has to be a $(1,0)$-identity  path, contradicting the nature of $\mathfrak{Q}$. 

Let us suppose, towards a contradiction, that there exists an index $k\in [1,m]$ for which 
$\mathfrak{R}_{k}$ is a  path of length strictly greater than one containing at least one echelon. Then, taking into account that $(\mathfrak{R}_{k},w_{k})\leq_{\mathbf{Pth}_{\boldsymbol{\mathcal{A}}}} (\mathfrak{P},s)$ and $\mathfrak{P}$ is an echelonless  path, we obtain, by Corollary~\ref{COrdII}, that 
$\bb{\mathfrak{R}_{k}}\leq\bb{\mathfrak{P}}$. Furthermore, since $(\mathfrak{Q},t)\leq_{\mathbf{Pth}_{\boldsymbol{\mathcal{A}}}} (\mathfrak{R}_{k},w_{k})$ and we are assuming that $\mathfrak{R}_{k}$ is a path of length strictly greater than one containing at least one  echelon, we obtain, according to Corollary~\ref{COrdI}, that $\bb{\mathfrak{Q}}<\bb{\mathfrak{R}_{k}}$. Finally, since $(\mathfrak{P},s)\leq_{\mathbf{Pth}_{\boldsymbol{\mathcal{A}}}} (\mathfrak{Q},t)$ and we are working under the assumption that $\mathfrak{Q}$ is an echelonless path, we obtain, by Corollary~\ref{COrdII}, that $\bb{\mathfrak{P}}\leq \bb{\mathfrak{Q}}$. This leads to the following contradiction
$$
\bb{\mathfrak{P}}\leq\bb{\mathfrak{Q}}<\bb{\mathfrak{R}_{k}}\leq \bb{\mathfrak{P}}.
$$
Therefore, we infer that, for every $k\in [1,m]$, $\mathfrak{R}_{k}$ cannot be a  path of length strictly greater than one containing at least one  echelon.

All in all, we conclude that for every $k\in [1,m]$, $\mathfrak{R}_{k}$ is an echelonless path. In particular $\mathfrak{P}$ is an echelonless path.

Now, since $(\mathfrak{P},s)<_{\mathbf{Pth}_{\boldsymbol{\mathcal{A}}}} (\mathfrak{Q},t)$, we have, by Remark~\ref{ROrd}, that there exists a natural number $n\in\mathbb{N}-\{0\}$, a word $\mathbf{v}\in S^{\star}$ of length $\bb{\mathbf{v}}=n+1$, and a family of paths $(\mathfrak{S}_{l})_{l\in\bb{\mathbf{v}}}$ in $\mathrm{Pth}_{\boldsymbol{\mathcal{A}},\mathbf{v}}$, such that $v_{0}=s$, $\mathfrak{S}_{0}=\mathfrak{P}$, $v_{n}=t$, $\mathfrak{S}_{n}=\mathfrak{P}$ and, for every $l\in n$, $(\mathfrak{S}_{l},v_{l})\prec_{\mathbf{Pth}_{\boldsymbol{\mathcal{A}}}} (\mathfrak{S}_{l+1},v_{l+1})$.

By means of a proof similar to the previous one, we can show that, for every $l\in [1,n]$, the  path $\mathfrak{S}_{l}$ must be an echelonless path. In particular, $\mathfrak{Q}$ is also an echelonless path.

Note that this situation fulfils the premises of Lemma~\ref{LOrdIII}. 

Assume that, for a unique word $\mathbf{c}\in S^{\star}-\{\lambda\}$, $\mathfrak{P}$ is a  $\mathbf{c}$-path in $\boldsymbol{\mathcal{A}}$  of the form
$$
\mathfrak{P}
=
\left(
(P_{i})_{i\in\bb{\mathbf{c}}+1},
(\mathfrak{p}_{i})_{i\in\bb{\mathbf{c}}},
(T_{i})_{i\in\bb{\mathbf{c}}}
\right),
$$
and, for a unique $\mathbf{d}\in S^{\star}$, $\mathfrak{Q}$ is a  $\mathbf{d}$-path in $\boldsymbol{\mathcal{A}}$  of the form
$$
\mathfrak{Q}
=
\left(
(Q_{j})_{j\in\bb{\mathbf{d}}+1},
(\mathfrak{q}_{j})_{j\in\bb{\mathbf{d}}},
(U_{j})_{j\in\bb{\mathbf{d}}}
\right).
$$ 

Since $(\mathfrak{Q},t)<_{\mathbf{Pth}_{\boldsymbol{\mathcal{A}}}} (\mathfrak{P},s)$ and the sequence witnessing this fact is entirely composed of echelonless  paths we have, by Lemma~\ref{LOrdIII}, that 
\[
\max\{
\bb{U_{j}}\mid j\in \bb{\mathbf{d}}
\}
<
\max\{
\bb{T_{i}}\mid i\in \bb{\mathbf{c}}
\}.
\]

On the other hand, since $(\mathfrak{P},s)<_{\mathbf{Pth}_{\boldsymbol{\mathcal{A}}}} (\mathfrak{Q},t)$ and the sequence witnessing this fact is entirely composed of  echelonless  paths we have, by Lemma~\ref{LOrdIII}, that 
\[
\max\{
\bb{T_{i}}\mid i\in \bb{\mathbf{c}}
\}
<
\max\{
\bb{U_{j}}\mid j\in \bb{\mathbf{d}}
\}
.
\]

This leads to the following contradiction
$$
\max\{
\bb{T_{i}}\mid i\in \bb{\mathbf{c}}
\}
<
\max\{
\bb{U_{j}}\mid j\in \bb{\mathbf{d}}
\}
<
\max\{
\bb{T_{i}}\mid i\in \bb{\mathbf{c}}
\}.
$$

Therefore, we can affirm that $\mathfrak{P}=\mathfrak{Q}$.

From this it follows that $\leq_{\coprod\mathrm{Pth_{\boldsymbol{\mathcal{A}}}}}$ is a partial order.

We next prove that the order $\leq_{\mathbf{Pth}_{\boldsymbol{\mathcal{A}}}}$ is Artinian, i.e., that there is not any strictly decreasing $\omega_{0}$-chain in the partially ordered set $(\coprod\mathrm{Pth}_{\boldsymbol{\mathcal{A}}}, \leq_{\mathbf{Pth}_{\boldsymbol{\mathcal{A}}}})$. Let us suppose, towards a contradiction, that there exists at least one strictly decreasing $\omega_{0}$-chain $((\mathfrak{R}_{k},s_{k}))_{k\in\omega_{0}}$ in $(\coprod\mathrm{Pth}_{\boldsymbol{\mathcal{A}}},\leq_{\mathbf{Pth}_{\boldsymbol{\mathcal{A}}}})$. Then, for every $k\in\omega_{0}$, we have that 
\begin{enumerate}
\item $\mathfrak{R}_{k}$ is a path in $\mathrm{Pth}_{\boldsymbol{\mathcal{A}},s_{k}}$, and
\item $(\mathfrak{R}_{k+1},s_{k+1})\prec_{\mathbf{Pth}_{\boldsymbol{\mathcal{A}}}}(\mathfrak{R}_{k},s_{k})$.
\end{enumerate}

Note that an strictly decreasing $\omega_{0}$-chain in $(\coprod\mathrm{Pth}_{\boldsymbol{\mathcal{A}}},\leq_{\mathbf{Pth}_{\boldsymbol{\mathcal{A}}}})$, in principle, should be defined by using, in clause (2) above, $<_{\mathbf{Pth}_{\boldsymbol{\mathcal{A}}}}$ instead of $\prec_{\mathbf{Pth}_{\boldsymbol{\mathcal{A}}}}$. However, by Remark~\ref{ROrd}, such a chain would, ultimately, lead to a chain as above, where $\prec_{\mathbf{Pth}_{\boldsymbol{\mathcal{A}}}}$ is actually used.

By construction, for every $k\in\omega_{0}$, the labelled  path $(\mathfrak{R}_{k},s_{k})$ cannot be minimal in the partially ordered set $(\coprod\mathrm{Pth}_{\boldsymbol{\mathcal{A}}},\leq_{\mathbf{Pth}_{\boldsymbol{\mathcal{A}}}})$. Otherwise, by the minimality of $(\mathfrak{R}_{k},s_{k})$, the labelled  path $(\mathfrak{R}_{k+1},s_{k+1})$ would be equal to $(\mathfrak{R}_{k},s_{k})$, contradicting the fact that $\prec_{\mathbf{Pth}_{\boldsymbol{\mathcal{A}}}}$ and $\Delta_{\coprod\mathrm{Pth}_{\boldsymbol{\mathcal{A}}}}$ are disjoint relations. Thus, by Lemma~\ref{LOrdI}, for every $k\in\omega_{0}$, we have that
\begin{enumerate}
\item $\mathfrak{R}_{k}$ is a $(1,0)$-identity path on a non-minimal term; or
\item $\mathfrak{R}_{k}$ is a  path of length strictly greater than one containing at least one echelon; or
\item $\mathfrak{R}_{k}$ is an echelonless path.
\end{enumerate}

We claim that in the sequence $((\mathfrak{R}_{k},s_{k}))_{k\in\omega_{0}}$, no $(1,0)$-identity  path can occur. Let us suppose, towards a contradiction, that this is not the case. Thus, assume that, for some $k\in\omega_{0}$, $\mathfrak{R}_{k}$ is a $(1,0)$-identity path, then the subsequence $((\mathfrak{R}_{l},s_{l}))_{l\geq k}$ is composed entirely of $(1,0)$-identity paths according to Proposition~\ref{PLower}. Hence, for every $l\in\omega_{0}$ with $l\geq k$, the  path $\mathfrak{R}_{l}$ has the form
$$
\mathfrak{R}_{l}=\mathrm{ip}^{(1,0)\sharp}_{s_{l}}\left(
R_{l}
\right),
$$
for a suitable  term $R_{l}$ in $\mathrm{T}_{\Sigma}(X)_{s_{l}}$.

Taking into account that we are assuming that, for every $l\in\omega_{0}$ with $l\geq k$,
$$
\left(\mathfrak{R}_{l+1},s_{l+1}
\right)
\prec_{\mathbf{Pth}_{\boldsymbol{\mathcal{A}}}}
\left(\mathfrak{R}_{l},s_{l}
\right)
$$
we have, according to Definition~\ref{DOrd} that, for every $l\in\omega_{0}$ with $l\geq k$, the following inequality holds
$$
\left(
R_{l+1},s_{l+1}
\right)
<_{
\mathbf{T}_{\Sigma}(X)
}
\left(
R_{l},s_{l}
\right).
$$
Thus $((R_{l},s_{l}))_{l\geq k}$ is an strictly decreasing $\omega_{0}$-chain in 
$
(\coprod\mathrm{T}_{\Sigma}(X),
\leq_{\mathbf{T}_{\Sigma}(X)})
$, 
contradicting the fact that $
(\coprod\mathrm{T}_{\Sigma}(X),
\leq_{\mathbf{T}_{\Sigma}(X)})
$ is an Artinian ordered set.

Thus, for every $k\in\omega_{0}$, $\mathfrak{R}_{k}$ cannot be a $(1,0)$-identity path.

We claim that, in the sequence $((\mathfrak{R}_{k},s_{k}))_{k\in\omega_{0}}$, the number of paths of length strictly greater than one containing at least one echelon is bounded. Let us suppose towards a contradiction that this is not the case. Then we can extract from $((\mathfrak{R}_{k},s_{k}))_{k\in\omega_{0}}$ an infinite subsequence $((\mathfrak{R}_{\varphi(k)},s_{\varphi(k)}))_{k\in\omega_{0}}$, where $\varphi$ is an strictly increasing endomapping of $\omega_{0}$, such that, for every $k\in\omega_{0}$, we have that
\begin{enumerate}
\item $\mathfrak{R}_{\varphi(k)}$ is a path of length strictly greater than one containing at least one echelon; and
\item $(\mathfrak{R}_{\varphi(k+1)}, s_{\varphi(k+1)})
<_{\mathbf{Pth}_{\boldsymbol{\mathcal{A}}}}
(\mathfrak{R}_{\varphi(k)}, s_{\varphi(k)})
$.
\end{enumerate}

But, by Corollary~\ref{COrdI}, this leads to a contradiction, concretely, to that of having an strictly decreasing $\omega_{0}$-chain in $(\mathbb{N},\leq)$ constructed from the length of such  paths.

From this, it follows that, in the sequence $((\mathfrak{R}_{k},s_{k}))_{k\in\omega_{0}}$, the number of indexes $k\in\omega_{0}$ for which $\mathfrak{R}_{k}$ is a  path of length strictly greater than one containing at least one  echelon is bounded.  

Let $l\in\omega_{0}$ be the maximum of the indexes for which $\mathfrak{R}_{l}$ is  a path of length strictly greater than one containing at least one  echelon. Then, for the infinite subsequence $((\mathfrak{R}_{k},s_{k}))_{k\in\omega_{0}-(l+1)}$ of $((\mathfrak{R}_{k},s_{k}))_{k\in\omega_{0}}$ and for every $k\in\omega_{0}-(l+1)$, we have that 
\begin{enumerate}
\item $\mathfrak{R}_{k}$ is an echelonless path; and
\item $(\mathfrak{R}_{k+1}, s_{k+1})
\prec_{\mathbf{Pth}_{\boldsymbol{\mathcal{A}}}}
(\mathfrak{R}_{k}, s_{k})
$.
\end{enumerate}
But by Lemma~\ref{LOrdIII}, this leads to a contradiction, concretely to that of having an strictly decreasing $\omega_{0}$-chain in $(\mathbb{N},\leq)$ constructed from the maximum of the heights of the sequence of translations of such paths.

Thus, no strictly decreasing $\omega_{0}$-chain can exist in $(\coprod\mathrm{Pth}_{\boldsymbol{\mathcal{A}}}, \leq_{\mathbf{Pth}_{\boldsymbol{\mathcal{A}}}})$. That is, the order $\leq_{\mathbf{Pth}_{\boldsymbol{\mathcal{A}}}}$ on $\coprod\mathrm{Pth}_{\boldsymbol{\mathcal{A}}}$ is Artinian.
\end{proof}

We next show that $\coprod\mathrm{ip}^{(1,0)\sharp}$ is an order embedding of the Artinian ordered set $(\coprod\mathrm{T}_{\Sigma}(X), \leq_{\mathbf{T}_{\Sigma}(X)})$, introduced on Remark~\ref{RTermOrd}, into the Artinian ordered set $(\coprod\mathrm{Pth}_{\boldsymbol{\mathcal{A}}}, \leq_{\mathbf{Pth}_{\boldsymbol{\mathcal{A}}}})$, introduced on Definition~\ref{DOrd}.

\begin{restatable}{proposition}{POrdEmb}
\label{POrdEmb} The mapping $\coprod\mathrm{ip}^{(1,0)\sharp}$ is an order embedding
$$
\textstyle
\coprod\mathrm{ip}^{(1,0)\sharp}\colon
\left(
\coprod\mathrm{T}_{\Sigma}(X)
, 
\leq_{\mathbf{T}_{\Sigma}(X)}
\right)
\mor 
\left(
\coprod\mathrm{Pth}_{\boldsymbol{\mathcal{A}}}, 
\leq_{\mathbf{Pth}_{\boldsymbol{\mathcal{A}}}}
\right).
$$
\end{restatable}
\begin{proof}
Let $s,t$ be sorts in $S$ and let $(Q,t)$ and let $(P,s)$ be pairs in $\coprod\mathrm{T}_{\Sigma}(X)$. We need to prove that the following statements are equivalent
\begin{enumerate}
\item $(Q,t)\leq_{\mathbf{T}_{\Sigma}(X)} (P,s)$;
\item $(\mathrm{ip}^{(1,0)\sharp}_{t}(Q),t)
\leq_{\mathbf{Pth}_{\boldsymbol{\mathcal{A}}}}
(\mathrm{ip}^{(1,0)\sharp}_{s}(P),s).
$ 
\end{enumerate}
Note that the statement trivially holds when the above statement are equalities. Therefore, it suffices to prove the above equivalence for the case of strict inequalities.

On the one hand, if we assume $(Q,t)<_{\mathbf{T}_{\Sigma}(X)} (P,s)$, then, in virtue of Definition~\ref{DOrd}, we have that $(\mathrm{ip}^{(1,0)\sharp}_{t}(Q),t)
\prec_{\mathbf{Pth}_{\boldsymbol{\mathcal{A}}}}
(\mathrm{ip}^{(1,0)\sharp}_{s}(P),s).
$ 

On the other hand, if we assume that $(\mathrm{ip}^{(1,0)\sharp}_{t}(Q),t)
\leq_{\mathbf{Pth}_{\boldsymbol{\mathcal{A}}}}
(\mathrm{ip}^{(1,0)\sharp}_{s}(P),s).
$ then, by Corollary~\ref{CLower}, we have that $(Q,t)<_{\mathbf{T}_{\Sigma}(X)} (P,s)$.

This finishes the proof.
\end{proof}

\chapter{
\texorpdfstring
{The categorial signature $\Sigma^{\boldsymbol{\mathcal{A}}}$ determined by $\boldsymbol{\mathcal{A}}$}
{The categorial signature}
}\label{S1D}

In this chapter, we define a new $S$-sorted signature, the categorial signature determined by 
$\boldsymbol{\mathcal{A}}$, denoted by $\Sigma^{\boldsymbol{\mathcal{A}}}$. This new signature extends the original signature $\Sigma$ by adding the operations of $0$-composition, $0$-source and $0$-target and, for every sort $s\in S$ and every rewrite rule $\mathfrak{p}\in \mathcal{A}_{s}$, a new constant operation symbol associated to it. This extension allows us to compare the free $S$-sorted $\Sigma$-algebra $\mathbf{T}_{\Sigma}(X)$ and the $\Sigma$-reduct of the free $S$-sorted $\Sigma^{\boldsymbol{\mathcal{A}}}$-algebra $\mathbf{T}_{\Sigma^{\boldsymbol{\mathcal{A}}}}(X)$. Moreover, we show that there exists a $\Sigma$-embedding $\eta^{(1,0)\sharp}$ of $\mathbf{T}_{\Sigma}(X)$ into $\mathbf{T}_{\Sigma^{\boldsymbol{\mathcal{A}}}}(X)$, which makes it possible to view all terms in $\mathbf{T}_{\Sigma}(X)$ as terms in $\mathbf{T}_{\Sigma^{\boldsymbol{\mathcal{A}}}}(X)$. In addition, the new categorical signature allows us to consider an $S$-sorted partial $\Sigma^{\boldsymbol{\mathcal{A}}}$-algebra structure on the set of many-sorted paths, denoted by $\mathbf{Pth}_{\boldsymbol{\mathcal{A}}}$, with the natural interpretations of the new operation symbols. We conclude this chapter by defining, by Artinian recursion, the $S$-sorted Curry-Howard mapping, denoted by $\mathrm{CH}^{(1)}$, from $\mathfrak{P}$ to $\mathbf{T}_{\Sigma^{\boldsymbol{\mathcal{A}}}}(X)_{s}$, which,  for every sort $s\in S$, sends a path $\mathfrak{P}$ in $\mathrm{Pth}_{\boldsymbol{\mathcal{A}},s}$ to a term 
$\mathrm{CH}^{(1)}_{s}(\mathfrak{P})$ in $\mathbf{T}_{\Sigma^{\boldsymbol{\mathcal{A}}}}(X)_{s}$. We then show that $\mathrm{CH}^{(1)}$ is a $\Sigma$-homomorphism but not a $\Sigma^{\boldsymbol{\mathcal{A}}}$-homomorphism and that $\coprod\mathrm{CH}^{(1)}$ is an order-preserving mapping.



We begin by defining the categorial signature determined  by $\mathcal{A}$ on $\Sigma$.

\begin{restatable}{definition}{DSigmaCat}
\label{DSigmaCat}\index{categorial signature!determined by $\mathcal{A}$, $\Sigma^{\boldsymbol{\mathcal{A}}}$}
Let $\Sigma$ be an $S$-sorted signature. Then the \emph{categorial signature determined by $\mathcal{A}$ on $\Sigma$}, denoted by $\Sigma^{\boldsymbol{\mathcal{A}}}$, is the $S$-sorted signature defined, for every $(\mathbf{s},s)\in S^{\star}\times S$, as follows:
$$
{\Sigma}^{\boldsymbol{\mathcal{A}}}_{\mathbf{s},s}=
\begin{cases}
{\Sigma}_{\mathbf{s},s}, 	
& \text{if $\mathbf{s}\not\in \{\lambda, (s),(s,s)\}$;}\\
{\Sigma}_{\mathbf{s},s}\amalg\mathcal{A}_{s}, 	
& \text{if $\mathbf{s}= \lambda$;}\\
{\Sigma}_{\mathbf{s},s}\amalg\{\mathrm{sc}^{0}_{s}, \mathrm{tg}^{0}_{s}\}, 	
& \text{if $\mathbf{s}= (s)$;}\\
{\Sigma}_{\mathbf{s},s}\amalg\{\circ^{0}_{s}\}, 	
& \text{if $\mathbf{s}= (s,s)$.}
\end{cases}
$$
That is, ${\Sigma}^{\boldsymbol{\mathcal{A}}}$ is the expansion of $\Sigma$ obtained by adding, for every sort $s\in S$,
(1) as many constants of coarity $s$ as there are rewrite rules in $\mathcal{A}_{s}$, (2) two unary operation symbols $\mathrm{sc}^{0}_{s}$ and $\mathrm{tg}^{0}_{s}$, both of them of arity $(s)$ and coarity $s$, which will be interpreted as total unary operations, and (3) a binary operation symbol $\circ^{0}_{s}$, of arity $(s,s)$ and coarity $s$, which will be interpreted as a binary (partial or total, depending on the case at hand) operation.
\end{restatable}

We next investigate how the free $S$-sorted $\Sigma^{\boldsymbol{\mathcal{A}}}$-algebra $\mathbf{T}_{\Sigma^{\boldsymbol{\mathcal{A}}}}(X)$ is related to previous $S$-sorted sets. But before doing that, we next define the embeddings of $X$, the $S$-sorted sets of variables, and  $\mathcal{A}$, the $S$-sorted set of rewrite rules, into the free $S$-sorted  $\Sigma^{\boldsymbol{\mathcal{A}}}$-algebra $\mathbf{T}_{\Sigma^{\boldsymbol{\mathcal{A}}}}(X)$, and the $\Sigma$-reduct of the free $\Sigma^{\boldsymbol{\mathcal{A}}}$-algebra $\mathbf{T}_{\Sigma^{\boldsymbol{\mathcal{A}}}}(X)$.

\begin{restatable}{definition}{DEta}
\label{DEta} Let $X$ be an $S$-sorted set and let $\mathbf{T}_{\Sigma^{\boldsymbol{\mathcal{A}}}}(X)$ be the free $S$-sorted $\Sigma^{\boldsymbol{\mathcal{A}}}$-algebra on $X$. We will denote by
\index{terms!first-order!$\mathbf{T}_{\Sigma^{\boldsymbol{\mathcal{A}}}}(X)$}
\begin{enumerate}
\item $\eta^{(1,X)}$\index{inclusion!first-order!$\eta^{(1,X)}$} the $S$-sorted mapping from $X$ to $\mathrm{T}_{\Sigma^{\boldsymbol{\mathcal{A}}}}(X)$ such that, for every sort $s\in S$, sends an element $x\in X_{s}$ to the variable $x\in\mathrm{T}_{\Sigma^{\boldsymbol{\mathcal{A}}}}(X)_{s}$.
\item $\eta^{(1,\mathcal{A})}$\index{inclusion!first-order!$\eta^{(1,\mathcal{A})}$} the $S$-sorted mapping from $\mathcal{A}$ to $\mathrm{T}_{\Sigma^{\boldsymbol{\mathcal{A}}}}(X)$ such that, for every sort $s\in S$, sends a rewrite rule $\mathfrak{p}\in\mathcal{A}_{s}$ to the constant $\mathfrak{p}^{\mathrm{T}_{\Sigma^{\boldsymbol{\mathcal{A}}}}(X)}$.
\end{enumerate}
The above $S$-sorted mappings are depicted in Figure~\ref{FPthEmb}.
\end{restatable}

\begin{figure}
\begin{tikzpicture}
[ACliment/.style={-{To [angle'=45, length=5.75pt, width=4pt, round]}},scale=.8]
\node[] (x) at (0,0) [] {$X$};
\node[] (a) at (0,-3) [] {$\mathcal{A}$};
\node[] (T) at (6,-3) [] {$\mathrm{T}_{\Sigma^{\boldsymbol{\mathcal{A}}}}(X)$};
\draw[ACliment, bend left=10]  (x) to node [above right] {$\eta^{(1,X)}$} (T);
\draw[ACliment]  (a) to node [below] {$\eta^{(1,\mathcal{A})}$} (T);
\end{tikzpicture}
\caption{Embeddings relative to $X$ and $\mathcal{A}$ at layer 1.}
\label{FPthEmb}
\end{figure}

We next define the $\Sigma$-reduct of the free $\Sigma^{\boldsymbol{\mathcal{A}}}$-algebra $\mathbf{T}_{\Sigma^{\boldsymbol{\mathcal{A}}}}(X)$.

\begin{restatable}{definition}{DRed}
\label{DRed}
Let $\mathrm{in}^{\Sigma,(1,0)}$ be the canonical embedding of $\Sigma$ into
$\Sigma^{\boldsymbol{\mathcal{A}}}$. Then, by Proposition~\ref{PFunSig}, for the morphism $\mathbf{in}^{\Sigma,(1,0)} = (\mathrm{id}^{S},\mathrm{in}^{\Sigma,(1,0)})$ from $(S,\Sigma)$ to
$(S,\Sigma^{\boldsymbol{\mathcal{A}}})$ and the free $\Sigma^{\boldsymbol{\mathcal{A}}}$-algebra $\mathbf{T}_{\Sigma^{\boldsymbol{\mathcal{A}}}}(X)$, we will denote by $\mathbf{T}_{\Sigma^{\boldsymbol{\mathcal{A}}}}^{(0,1)}(X)$ the $\Sigma$-algebra $\mathbf{in}^{\Sigma,(0,1)}(\mathbf{T}_{\Sigma^{\boldsymbol{\mathcal{A}}}}(X))$. We will call $\mathbf{T}_{\Sigma^{\boldsymbol{\mathcal{A}}}}^{(0,1)}(X)$ the \emph{$\Sigma$-reduct} of the free $\Sigma^{\boldsymbol{\mathcal{A}}}$-algebra $\mathbf{T}_{\Sigma^{\boldsymbol{\mathcal{A}}}}(X)$. 
\end{restatable}

\begin{remark}\label{RRed} The underlying $S$-sorted of $\mathbf{T}_{\Sigma^{\boldsymbol{\mathcal{A}}}}^{(0,1)}(X)$ is the same as that of the free $S$-sorted $\Sigma^{\boldsymbol{\mathcal{A}}}$-algebra 
$\mathbf{T}_{\Sigma^{\boldsymbol{\mathcal{A}}}}(X)$, while, for every pair $(\mathbf{s},s)\in S^{\star}\times S$ and every operation symbol $\sigma\in\Sigma_{\mathbf{s},s}$,
$$
\sigma^{\mathbf{T}^{(0,1)}_{\Sigma^{\boldsymbol{\mathcal{A}}}}(X)}=\sigma^{\mathbf{T}_{\Sigma^{\boldsymbol{\mathcal{A}}}}(X)}.
$$
\end{remark}

\begin{restatable}{proposition}{PEmb}
\label{PEmb}
\index{inclusion!first-order!$\eta^{(1,0)\sharp}$}
There exists an embedding from the free $S$-sorted $\Sigma$-algebra $\mathbf{T}_{\Sigma}(X)$ into the 
$\Sigma$-reduct $\mathbf{T}^{(0,1)}_{\Sigma^{\boldsymbol{\mathcal{A}}}}(X)$ of $\mathbf{T}_{\Sigma^{\boldsymbol{\mathcal{A}}}}(X)$.
\end{restatable}
\begin{proof}
Consider $\eta^{(0,X)}$ and $\eta^{(1,X)}$, the canonical embeddings of $X$ into $\mathrm{T}_{\Sigma}(X)$ and $\mathrm{T}_{\Sigma^{\boldsymbol{\mathcal{A}}}}(X)$, respectively. Since $\mathbf{T}^{(0,1)}_{\Sigma^{\boldsymbol{\mathcal{A}}}}(X)$ is a $\Sigma$-algebra, by the universal property of the free $\Sigma$-algebra on $X$, i.e., $\mathbf{T}_{\Sigma}(X)$, there exists a unique many-sorted $\Sigma$-homomorphism, that we will denote by $\eta^{(1,0)\sharp}$ from $\mathbf{T}_{\Sigma}(X)$ to $\mathbf{T}^{(0,1)}_{\Sigma^{\boldsymbol{\mathcal{A}}}}(X)$ such that 
$$
\eta^{(1,0)\sharp}\circ\eta^{(0,X)}=\eta^{(1,X)}.
$$

The many-sorted $\Sigma$-homomorphism $\eta^{(1,0)\sharp}$ sends, for every sort $s\in S$, a term $P\in\mathrm{T}_{\Sigma}(X)_{s}$ to the term $P\in \mathrm{T}_{\Sigma^{\boldsymbol{\mathcal{A}}}}(X)_{s}$. This many-sorted mapping is injective. The reader is advised to consult the diagram presented in Figure~\ref{FEmb}.

As usual, we identify $\mathrm{T}_{\Sigma}(X)$ and $\mathrm{Im}[\eta^{(1,0)\sharp}]$. In this way, $\mathrm{T}_{\Sigma}(X)$ becomes a subset of $\mathrm{T}_{\Sigma^{\boldsymbol{\mathcal{A}}}}(X)$ and $\mathbf{T}_{\Sigma}(X)$ a subalgebra of $\mathbf{T}_{\Sigma^{\boldsymbol{\mathcal{A}}}}^{(0,1)}(X)$.
\end{proof}

\begin{figure}
\begin{center}
\begin{tikzpicture}
[ACliment/.style={-{To [angle'=45, length=5.75pt, width=4pt, round]}
}, scale=0.8]
\node[] (xq) 		at 	(0,0) 	[] 	{$X$};
\node[] (txq) 	at 	(6,0) 	[] 	{$\mathrm{T}_{\Sigma}(X)$};
\node[] (txqc) 	at 	(6,-3) 	[] 	{$\mathrm{T}^{(0,1)}_{\Sigma^{\boldsymbol{\mathcal{A}}}}
(X)$};
\draw[ACliment]  (xq) 	to node [above]	
{$\eta^{(0,X)}$} (txq);
\draw[ACliment]  (txq) 	to node [right]	
{$\eta^{(1,0)\sharp}$} (txqc);
\draw[ACliment, bend right=10] (xq) 	to node [below left]	
{$\eta^{(1,X)}$} (txqc);
\end{tikzpicture}
\end{center}
\caption{Embeddings relative to $X$ at layers 0 \& 1.}
\label{FEmb}
\end{figure}

\section{
\texorpdfstring
{A structure of partial $\Sigma^{\boldsymbol{\mathcal{A}}}$-algebra on $\mathrm{Pth}_{\boldsymbol{\mathcal{A}}}$}
{A partial algebra on paths}
}
We next show that the $S$-sorted set $\mathrm{Pth}_{\boldsymbol{\mathcal{A}}}$ has a natural structure of 
$S$-sorted partial $\Sigma^{\boldsymbol{\mathcal{A}}}$-algebra.

\begin{restatable}{proposition}{PPthCatAlg}
\label{PPthCatAlg}
\index{path!first-order!$\mathbf{Pth}_{\boldsymbol{\mathcal{A}}}$}
The $S$-sorted set $\mathrm{Pth}_{\boldsymbol{\mathcal{A}}}$ is equipped, in a natural way, with a structure of 
$S$-sorted partial $\Sigma^{\boldsymbol{\mathcal{A}}}$-algebra.
\end{restatable}

\begin{proof}
Let us denote by $\mathbf{Pth}_{\boldsymbol{\mathcal{A}}}$ the $\Sigma^{\boldsymbol{\mathcal{A}}}$-algebra defined as follows:

\textsf{(1)} The underlying $S$-sorted set of $\mathbf{Pth}_{\boldsymbol{\mathcal{A}}}$ is $\mathrm{Pth}_{\boldsymbol{\mathcal{A}}} = (\mathrm{Pth}_{\boldsymbol{\mathcal{A}},s})_{s\in S}$.

\textsf{(2)} For every $(\mathbf{s},s)\in S^{\star}\times S$ and every operation symbol $\sigma\in\Sigma_{\mathbf{s},s}$, the operation $\sigma^{\mathbf{Pth}_{\boldsymbol{\mathcal{A}}}}$ is given by the interpretation of $\sigma$ in $\mathbf{Pth}^{(0,1)}_{\boldsymbol{\mathcal{A}}}$ that, we recall, was stated in Proposition~\ref{PPthAlg}, where, in addition, we proved that $\sigma^{\mathbf{Pth}_{\boldsymbol{\mathcal{A}}}}$ is well-defined. 

\textsf{(3)} For every $s\in S$ and every $\mathfrak{p}\in\mathcal{A}_{s}$ with $\mathfrak{p}=(M,N)$, the constant $\mathfrak{p}^{\mathbf{Pth}_{\boldsymbol{\mathcal{A}}}}$ is given by the echelon determined by $\mathfrak{p}$, i.e., by the one-step path of sort $s$ $\mathrm{ech}^{(1,\mathcal{A})}_{s}(\mathfrak{p})$ that has $M$ as $(0,1)$-source and, after applying the rewrite rule $\mathfrak{p}$ retrieves $N$ as $(0,1)$-target.  That is, the constant $\mathfrak{p}$ is interpreted as the following echelon:
\begin{center}
\begin{tikzpicture}
[ACliment/.style={-{To [angle'=45, length=5.75pt, width=4pt, round]},
font=\scriptsize}]
\node[] (1) at (0,0) [] {$\mathrm{ech}^{(1,\mathcal{A})}_{s}(\mathfrak{p})\colon M$};
\node[] (2) at (3.5,0) [] {$N$.};
\draw[ACliment]  (1) to node [above]
{$(\mathfrak{p},\mathrm{id}^{\mathrm{T}_{\Sigma}(X)_{s}})$} (2);
\end{tikzpicture}
\end{center}

\textsf{(4)} For every $s\in S$, the operation $\mathrm{sc}_{s}^{0\mathbf{Pth}_{\boldsymbol{\mathcal{A}}}}$ is equal to $\mathrm{ip}^{(1,0)\sharp}_{s}\circ\mathrm{sc}^{(0,1)}_{s}$, i.e., to the mapping  that assigns to a path $\mathfrak{P}$ in $\mathrm{Pth}_{\boldsymbol{\mathcal{A}}, s}$ the $(1,0)$-identity path on the $(0,1)$-source of $\mathfrak{P}$.

\textsf{(5)} For every $s\in S$, the operation $\mathrm{tg}_{s}^{0\mathbf{Pth}_{\boldsymbol{\mathcal{A}}}}$ is equal to $\mathrm{ip}^{(1,0)\sharp}_{s}\circ\mathrm{tg}^{(0,1)}_{s}$, i.e., to the mapping  that assigns to a path $\mathfrak{P}$ in $\mathrm{Pth}_{\boldsymbol{\mathcal{A}}, s}$  the $(1,0)$-identity path on the $(0,1)$-target of $\mathfrak{P}$.  

\textsf{(6)} For every $s\in S$, the partial binary operation $\circ_{s}^{0\mathbf{Pth}_{\boldsymbol{\mathcal{A}}}}$ is equal to the partial binary operation of $0$-composition of paths that, we recall, was introduced in Definition~\ref{DPthComp}.  We recall that, in Proposition~\ref{PPthComp}, we proved that given two paths $\mathfrak{P}$ and $\mathfrak{Q}$ such that 
$$
\mathrm{sc}^{(0,1)}_{s}\left(\mathfrak{Q}\right)=
\mathrm{tg}^{(0,1)}_{s}\left(\mathfrak{P}\right),
$$
then the $0$-composition  $\mathfrak{Q}\circ_{s}^{0\mathbf{Pth}_{\boldsymbol{\mathcal{A}}}}\mathfrak{P}$ is well-defined.

This completes the definition of the partial $\Sigma^{\boldsymbol{\mathcal{A}}}$-algebra $\mathbf{Pth}_{\boldsymbol{\mathcal{A}}}$.
\end{proof}

\begin{remark}\label{RPthCatAlg} For the partial $\Sigma^{\boldsymbol{\mathcal{A}}}$-algebra $\mathbf{Pth}_{\boldsymbol{\mathcal{A}}}$, we denote by $\mathbf{Pth}^{(0,1)}_{\boldsymbol{\mathcal{A}}}$ the $\Sigma$-algebra $\mathbf{in}^{\Sigma,(0,1)}(\mathbf{Pth}_{\boldsymbol{\mathcal{A}}})$. We will call $\mathbf{Pth}^{(0,1)}_{\boldsymbol{\mathcal{A}}}$ the \emph{$\Sigma$-reduct} of the partial $\Sigma^{\boldsymbol{\mathcal{A}}}$-algebra $\mathbf{Pth}_{\boldsymbol{\mathcal{A}}}$. Note that this $\Sigma$-algebra coincides with the $\Sigma$-algebra introduced in Proposition~\ref{PPthAlg}.
\end{remark}

\section{The Curry-Howard mapping}

The previous results, will allow us, by means of an $S$-sorted mapping, to consider paths in the rewriting  system $\boldsymbol{\mathcal{A}} =(S,\Sigma,X,\mathcal{A})$ as terms relative to $\Sigma^{\boldsymbol{\mathcal{A}}}$ and $X$. To do this, we will define, by Artinian recursion, an $S$-sorted mapping from $\mathrm{Pth}_{\boldsymbol{\mathcal{A}}}$ to $\mathrm{T}_{\Sigma^{\boldsymbol{\mathcal{A}}}}(X)$, the underlying $S$-sorted set of $\mathbf{T}_{\Sigma^{\boldsymbol{\mathcal{A}}}}(X)$, the free $\Sigma^{\boldsymbol{\mathcal{A}}}$-algebra on $X$. In this way, every path in $\boldsymbol{\mathcal{A}}$ will be denoted by a term in $\mathrm{T}_{\Sigma^{\boldsymbol{\mathcal{A}}}}(X)$. Since this mapping reminds us of the classical Curry-Howard correspondence (see~\cite{CF58} and \cite{H80}), we have decided to denote it by $\mathrm{CH}^{(1)}$.

\begin{restatable}{definition}{DCH}
\label{DCH}  
\index{Curry-Howard!first-order!$\mathrm{CH}^{(1)}$}
The \emph{Curry-Howard mapping} is the $S$-sorted mapping
$$
\textstyle
\mathrm{CH}^{(1)}\colon\mathrm{Pth}_{\boldsymbol{\mathcal{A}}}\mor\mathrm{T}_{\Sigma^{\boldsymbol{\mathcal{A}}}}(X)
$$
defined by Artinian recursion on $(\coprod\mathrm{Pth}_{\boldsymbol{\mathcal{A}}}, \leq_{\mathbf{Pth}_{\boldsymbol{\mathcal{A}}}})$ as follows.

\textsf{Base step of the Artinian recursion}.

Let $(\mathfrak{P},s)$ be a minimal element of $(\coprod\mathrm{Pth}_{\boldsymbol{\mathcal{A}}}, \leq_{\mathbf{Pth}_{\boldsymbol{\mathcal{A}}}})$. Then, by Proposition~\ref{PMinimal}, the path $\mathfrak{P}$ is either~(1) an $(1,0)$-identity path or~(2) an echelon.

If~(1), i.e., if $\mathfrak{P}$ is an $(1,0)$-identity path, then $\mathfrak{P}=\mathrm{ip}^{(1,0)\sharp}_{s}(P)$ for some term $P\in\mathrm{T}_{\Sigma}(X)_{s}$. We define $\mathrm{CH}^{(1)}_{s}(\mathfrak{P})$ to be the term in $\mathrm{T}_{\Sigma^{\boldsymbol{\mathcal{A}}}}(X)_{s}$ given by the lift of the term $P$ by means of the embedding $\eta^{(1,0)\sharp}$, i.e.,
$$
\mathrm{CH}^{(1)}_{s}
\left(
\mathfrak{P}
\right)=\eta^{(1,0)\sharp}_{s}(P).
$$

If~(2), i.e., if $\mathfrak{P}$ is an echelon associated to a rewrite rule $\mathfrak{p}=(M,N)$, that is if $\mathfrak{P}$ has the form
$$
\xymatrix@C=55pt{
\mathfrak{P}: M
\ar[r]^-{\text{\Small{$(\mathfrak{p},\mathrm{id}^{\mathrm{T}_{\Sigma}(X)_{s}})$}}}
&
N
},
$$
then we define $\mathrm{CH}^{(1)}_{s}(\mathfrak{P})$ to be the syntactic representation of the unique rewrite rule occurring in $\mathfrak{P}$, i.e., 
$$
\mathrm{CH}^{(1)}_{s}\left(\mathfrak{P}
\right)=
\mathfrak{p}^{\mathbf{T}_{\Sigma^{\boldsymbol{\mathcal{A}}}}(X)}.
$$

\textsf{Inductive step of the Artinian recursion}.

Let $(\mathfrak{P},s)$ be a non-minimal element of $(\coprod\mathrm{Pth}_{\boldsymbol{\mathcal{A}}}, \leq_{\mathbf{Pth}_{\boldsymbol{\mathcal{A}}}})$. We can assume that $\mathfrak{P}$ is a not a $(1,0)$-identity path, since those paths already have an image for the Curry-Howard mapping. Let us suppose that, for every sort $t\in S$ and every path $\mathfrak{Q}\in\mathrm{Pth}_{\boldsymbol{\mathcal{A}},t}$, if $(\mathfrak{Q},t)<_{\mathbf{Pth}_{\boldsymbol{\mathcal{A}}}}(\mathfrak{P},s)$, then the value of the Curry-Howard mapping at $\mathfrak{Q}$, i.e., $\mathrm{CH}^{(1)}_{t}(\mathfrak{Q})$, has already been defined.

By Lemma~\ref{LOrdI}, we have that $\mathfrak{P}$ is either~(1) a path of length strictly greater than one containing at least one echelon or~(2) an echelonless path.

If~(1), i.e., if $\mathfrak{P}$ is a path of length strictly greater than one containing at least on echelon, then let $i\in \bb{\mathfrak{P}}$ be the first index for which the one-step subpath $\mathfrak{P}^{i,i}$ of $\mathfrak{P}$ is an echelon. We consider different cases for $i$ according to the cases presented in Definition~\ref{DOrd}.

If $i=0$, we have that the pairs $(\mathfrak{P}^{0,0},s)$ and $(\mathfrak{P}^{1,\bb{\mathfrak{P}}-1},s)$ $\prec_{\mathbf{Pth}_{\boldsymbol{\mathcal{A}}}}$-precede the pair $(\mathfrak{P},s)$. Therefore, the values of the Curry-Howard mapping at $\mathfrak{P}^{0,0}$ and $\mathfrak{P}^{1,\bb{\mathfrak{P}}-1}$, respectively, have already been defined. 

In this case, we set
$$
\mathrm{CH}^{(1)}_{s}\left(
\mathfrak{P}
\right)=
\mathrm{CH}^{(1)}_{s}\left(
\mathfrak{P}^{1,\bb{\mathfrak{P}}-1}
\right)
\circ_{s}^{0\mathbf{T}_{\Sigma^{\boldsymbol{\mathcal{A}}}}(X)}
\mathrm{CH}^{(1)}_{s}\left(
\mathfrak{P}^{0,0}
\right).
$$

If $i\neq 0$, we have that the pairs $(\mathfrak{P}^{0,i-1},s)$ and $(\mathfrak{P}^{i,\bb{\mathfrak{P}}-1},s)$ $\prec_{\mathbf{Pth}_{\boldsymbol{\mathcal{A}}}}$-precede the pair $(\mathfrak{P},s)$. Therefore, the values of the Curry-Howard mapping at $\mathfrak{P}^{0,i-1}$ and $\mathfrak{P}^{i,\bb{\mathfrak{P}}-1}$, respectively, have already been defined. 

In this case, we set
$$
\mathrm{CH}^{(1)}_{s}\left(
\mathfrak{P}
\right)=
\mathrm{CH}^{(1)}_{s}\left(
\mathfrak{P}^{i,\bb{\mathfrak{P}}-1}
\right)
\circ_{s}^{0\mathbf{T}_{\Sigma^{\boldsymbol{\mathcal{A}}}}(X)}
\mathrm{CH}^{(1)}_{s}\left(
\mathfrak{P}^{0,i-1}
\right).
$$

This finishes the definition of the value of the Curry-Howard mapping at a path of length strictly greater than one containing at least one echelon.

If~(2), i.e., if $\mathfrak{P}$ is an echelonless path in $\mathrm{Pth}_{\boldsymbol{\mathcal{A}},s}$, then the conditions for the path extraction algorithm, as stated in Lemma~\ref{LPthExtract}, are fulfilled. Then, by Lemma~\ref{LPthHeadCt}, there exists a unique word $\mathbf{s}\in S^{\star}-\{\lambda\}$ and a unique operation symbol $\sigma\in \Sigma_{\mathbf{s},s}$ associated to $\mathfrak{P}$. Let $(\mathfrak{P}_{j})_{j\in\bb{\mathbf{s}}}$ be the family of paths in $\mathrm{Pth}_{\boldsymbol{\mathcal{A}},\mathbf{s}}$ which, in virtue of Lemma~\ref{LPthExtract}, we can extract from $\mathfrak{P}$. Note that, for every $j\in\bb{\mathbf{s}}$, we have that 
$(\mathfrak{P}_{j},s_{j})$ $\prec_{\mathbf{Pth}_{\boldsymbol{\mathcal{A}}}}$ $(\mathfrak{P},s)$.  Therefore, for every $j\in\bb{\mathbf{s}}$, the value of the Curry-Howard mapping at $\mathfrak{P}_{j}$ has already been defined.

In this case, we set
$$
\mathrm{CH}^{(1)}_{s}\left(\mathfrak{P}
\right)=
\sigma^{\mathbf{T}_{\Sigma^{\boldsymbol{\mathcal{A}}}}(X)}
\left(\left(\mathrm{CH}^{(1)}_{s_{j}}\left(
\mathfrak{P}_{j}
\right)\right)_{j\in\bb{\mathbf{s}}}\right).
$$

That finishes the definition of the value of the Curry-Howard mapping at an echelonless path.

This completes the definition of the Curry-Howard mapping.
\end{restatable}

\begin{remark}\label{RCH} The Curry-Howard mapping associates to each path of sort $s\in S$ a term in 
$\mathrm{T}_{\Sigma^{\boldsymbol{\mathcal{A}}}}(X)_{s}$. In this regard, we emphasize how crucial the categorial expansion of the original signature has been. The reader will also notice that the Curry-Howard mapping matches a path over terms with a term, but, at the cost of adding complexity on the term side. This term contains the most important features of the path and will be used later to reconstruct another path, not necessarily the same as the original but normalized, in which the rewrite rules acting in parallel are applied following a leftmost innermost derivation strategy.
\end{remark}

\section{
\texorpdfstring
{The behaviour of the Curry-Howard mapping}
{Behaviour}
}
The following propositions will provide a deeper understanding of the just defined Curry-Howard mapping.  Moreover, we study how the Curry-Howard mapping relates to the mappings from $X$ and $\mathcal{A}$.

We next prove that the Curry-Howard mapping is a $\Sigma$-homomorphism from the partial $\Sigma$-algebra $\mathbf{Pth}^{(0,1)}_{\boldsymbol{\mathcal{A}}}$ to the $\Sigma$-algebra $\mathbf{T}^{(0,1)}_{\Sigma^{\boldsymbol{\mathcal{A}}}}(X)$.

\begin{restatable}{proposition}{PCHHom}
\label{PCHHom}
The Curry-Howard mapping is a $\Sigma$-homomorphism from $\mathbf{Pth}^{(0,1)}_{\boldsymbol{\mathcal{A}}}$ to $\mathbf{T}^{(0,1)}_{\Sigma^{\boldsymbol{\mathcal{A}}}}(X)$. 
\end{restatable}
\begin{proof}
Let $(\mathbf{s},s)$ be an element of $S^{\star}\times S$ and $\sigma$ an operation symbol in $\Sigma_{\mathbf{s},s}$.

If $\mathbf{s}=\lambda$, then the following chain of equalities holds
\allowdisplaybreaks
\begin{align*}
\mathrm{CH}^{(1)}_{s}\left(
\sigma^{\mathbf{Pth}_{\boldsymbol{\mathcal{A}}}}
\right)
&=
\mathrm{CH}^{(1)}_{s}\left(
\mathrm{ip}^{(1,0)\sharp}_{s}\left(
\sigma^{\mathbf{T}_{\Sigma}(X)}
\right)
\right)
\tag{1}
\\&=
\eta^{(1,0)\sharp}_{s}\left(\sigma^{\mathbf{T}_{\Sigma}(X)}
\right)
\tag{2}
\\&=
\sigma^{\mathbf{T}_{\Sigma^{\boldsymbol{\mathcal{A}}}}(X)}.
\tag{3}
\end{align*}

The first equality follows from the fact that, for a constant symbol $\sigma$ in $\Sigma_{\lambda,s}$ then, following Remark~\ref{RConsSigma}, the interpretation of $\sigma$ in $\mathbf{Pth}_{\boldsymbol{\mathcal{A}}}$ is given by $
\sigma^{\mathbf{Pth}_{\boldsymbol{\mathcal{A}}}}=\mathrm{ip}^{(1,0)\sharp}_{s}(\sigma^{\mathbf{T}_{\Sigma}(X)})$; the second equality follows from Proposition~\ref{PCHId}; finally, the last equality follows from the fact that $\eta^{(1,0)\sharp}$ is a $\Sigma$-homomorphism.

We now consider the case in which $\mathbf{s}\neq\lambda$. Let $(\mathfrak{P}_{j})_{j\in\bb{\mathbf{s}}}$ be a family of paths in $\mathrm{Pth}_{\boldsymbol{\mathcal{A}},\mathbf{s}}$. We consider different cases according to the nature of the family $(\mathfrak{P}_{j})_{j\in\bb{\mathbf{s}}}$. It could be the case that either~(1), for every $j\in\bb{\mathbf{s}}$, $\mathfrak{P}_{j}$ is a $(1,0)$-identity path or~(2), there exists an index $j\in\bb{\mathbf{s}}$ for which $\mathfrak{P}_{j}$ is a non-$(1,0)$-identity path.

If~(1), then for every $j\in\bb{\mathbf{s}}$, $\mathfrak{P}_{j}$ is equal to $\mathrm{ip}^{(1,0)\sharp}_{s_{j}}(P_{j})$, for some term $P_{j}$ in $\mathrm{T}_{\Sigma}(X)_{s_{j}}$. In this case, the following chain of equalities holds
\allowdisplaybreaks
\begin{align*}
\mathrm{CH}^{(1)}_{s}\left(
\sigma^{\mathbf{Pth}_{\boldsymbol{\mathcal{A}}}}\left(
\left(\mathfrak{P}_{j}
\right)_{j\in\bb{\mathbf{s}}}\right)
\right)
&=
\mathrm{CH}^{(1)}_{s}\left(
\sigma^{\mathbf{Pth}_{\boldsymbol{\mathcal{A}}}}\left(
\left(
\mathrm{ip}^{(1,0)\sharp}_{s_{j}}\left(
P_{j}
\right)
\right)_{j\in\bb{\mathbf{s}}}\right)
\right)
\tag{1}
\\&=
\mathrm{CH}^{(1)}_{s}\left(
\mathrm{ip}^{(1,0)\sharp}_{s}\left(
\sigma^{\mathbf{T}_{\Sigma}(X)}\left(\left(
P_{j}
\right)_{j\in\bb{\mathbf{s}}}\right)
\right)\right)
\tag{2}
\\&=
\eta^{(1,0)\sharp}_{s}\left(
\sigma^{\mathbf{T}_{\Sigma}(X)}\left(\left(
P_{j}
\right)_{j\in\bb{\mathbf{s}}}\right)
\right)
\tag{3}
\\&=
\sigma^{\mathbf{T}_{\Sigma^{\boldsymbol{\mathcal{A}}}}(X)}\left(
\left(
\eta^{(1,0)\sharp}_{s_{j}}\left(
P_{j}
\right)\right)_{j\in\bb{\mathbf{s}}}\right)
\tag{4}
\\&=
\sigma^{\mathbf{T}_{\Sigma^{\boldsymbol{\mathcal{A}}}}(X)}\left(
\left(
\mathrm{CH}^{(1)}_{s_{j}}\left(
\mathrm{ip}^{(1,0)\sharp}_{s_{j}}\left(
P_{j}
\right)
\right)
\right)_{j\in\bb{\mathbf{s}}}\right)
\tag{5}
\\&=
\sigma^{\mathbf{T}_{\Sigma^{\boldsymbol{\mathcal{A}}}}(X)}\left(
\left(
\mathrm{CH}^{(1)}_{s_{j}}\left(
\mathfrak{P}_{j}
\right)
\right)_{j\in\bb{\mathbf{s}}}\right).
\tag{6}
\end{align*}

The first equality follows from the fact that, for every $j\in\bb{\mathbf{s}}$, $\mathfrak{P}_{j}$ is a $(1,0)$-identity path;  the second equality follows from the fact that, by Proposition~\ref{PIpHom} $\mathrm{ip}^{(1,0)\sharp}$ is a $\Sigma$-homomorphism from $\mathbf{T}_{\Sigma}(X)$ to $\mathbf{Pth}^{(0,1)}_{\boldsymbol{\mathcal{A}}}$;  the third equality follows from the characterization of the value of the Curry-Howard mapping on $(1,0)$-identity paths introduced in Proposition~\ref{PCHId}; the fourth equality follows from the fact that, by Proposition~\ref{PEmb}, $\eta^{(1,0)\sharp}$ is a $\Sigma$-homomorphism from $\mathbf{T}_{\Sigma}(X)$ to $\mathbf{T}^{(0,1)}_{\Sigma^{\boldsymbol{\mathcal{A}}}}(X)$; the fifth equality follows from the characterization of the value of the Curry-Howard mapping on $(1,0)$-identity paths introduced in Proposition~\ref{PCHId}; finally, the last equality recovers, for every $j\in\bb{\mathbf{s}}$, the definition of $\mathfrak{P}_{j}$ as a $(1,0)$-identity path.

This proves Case~(1).

If~(2), i.e., if there exists some index $j\in\bb{\mathbf{s}}$ for which $\mathfrak{P}_{j}$ is a non-$(1,0)$-identity  path then, according to Corollary~\ref{CPthWB}, $\sigma^{\mathbf{Pth}_{\boldsymbol{\mathcal{A}}}}((\mathfrak{P}_{j})_{j\in\bb{\mathbf{s}}})$ is an echelonless path. Moreover, according to Proposition~\ref{PRecov}, the path extraction algorithm from Lemma~\ref{LPthExtract} applied to it retrieves the original family $(\mathfrak{P}_{j})_{j\in\bb{\mathbf{s}}}$. Then, following Definition~\ref{DCH},  the Curry-Howard mapping applied to $\sigma^{\mathbf{Pth}_{\boldsymbol{\mathcal{A}}}}((\mathfrak{P}_{j})_{j\in\bb{\mathbf{s}}})$ is given by
\begin{align*}
\mathrm{CH}^{(1)}_{s}\left(
\sigma^{\mathbf{Pth}_{\boldsymbol{\mathcal{A}}}}
\left(\left(
\mathfrak{P}_{j}
\right)_{j\in\bb{\mathbf{s}}}
\right)\right)
&=
\sigma^{\mathbf{T}_{\Sigma^{\boldsymbol{\mathcal{A}}}}(X)}
\left(\left(
\mathrm{CH}^{(1)}_{s_{j}}
\left(
\mathfrak{P}_{j}
\right)\right)_{j\in\bb{\mathbf{s}}}
\right).
\end{align*}

Case~(2) follows.

This finishes the proof.
\end{proof}

\begin{figure}
\begin{center}
\begin{tikzpicture}
[ACliment/.style={-{To [angle'=45, length=5.75pt, width=4pt, round]}},scale=0.8]
\node[] (x) at (0,0) [] {$X$};
\node[] (t) at (6,0) [] {$\mathrm{T}_{\Sigma}(X)$};
\node[] (p) at (6,-2) [] {$\mathrm{Pth}_{\boldsymbol{\mathcal{A}}}$};
\node[] (tc) at (6,-4) [] {$\mathrm{T}_{\Sigma^{\boldsymbol{\mathcal{A}}}}(X)$};

\draw[ACliment]  (x) to node [above right] {$\eta^{(0,X)}$} (t);
\draw[ACliment, bend right=10]  (x) to node [midway, fill=white] {$\mathrm{ip}^{(1,X)}$} (p);
\draw[ACliment, bend right=20]  (x) to node [below left] {$\eta^{(1,X)}$} (tc);

\draw[ACliment]  (t) to node [right] {$\mathrm{ip}^{(1,0)\sharp}$} (p);
\draw[ACliment]  (p) to node [right] {$\mathrm{CH}^{(1)}$} (tc);
\draw[ACliment, rounded corners] (t.east)
--
 ($(tc.east)+(.5,4)$)
--  node [right] {$\eta^{(1,0)\sharp}$} 
($(tc.east)+(.5,0)$)
-- (tc.east);
\end{tikzpicture}
\end{center}
\caption{Behaviour of $\mathrm{CH}^{(1)}$ relative to $X$ at layers 0 \& 1.}
\label{FCHId}
\end{figure}

In the following proposition we describes the image under the Curry-Howard mapping of $(1,0)$-identity paths.

\begin{proposition}\label{PCHId}
For the Curry-Howard mapping we have that
\begin{multicols}{2}
\begin{itemize}
\item[(i)] $\mathrm{CH}^{(1)}\circ\mathrm{ip}^{(1,X)}=\eta^{(1,X)}$,
\item[(ii)] $\mathrm{CH}^{(1)}\circ\mathrm{ip}^{(1,0)\sharp}=\eta^{(1,0)\sharp}$.
\end{itemize}
\end{multicols}
\end{proposition}
\begin{proof}
The reader is advised to consult the diagram presented in Figure~\ref{FCHId}.

Let $s$ be a sort in $S$ and let $x$ be a variable in $X_{s}$. Let us note that the following chain of equalities holds
\allowdisplaybreaks
\begin{align*}
\mathrm{CH}^{(1)}_{s}\left(
\mathrm{ip}^{(1,X)}_{s}(x)
\right)&=
\mathrm{CH}^{(1)}_{s}\left(x,\lambda,\lambda
\right)
\tag{1}
\\&=\eta^{(1,0)\sharp}_{s}(x)
\tag{2}
\\&=\eta^{(1,0)\sharp}_{s}\left(
\eta^{(0,X)}_{s}(x)
\right)
\tag{3}
\\&=\eta^{(1,X)}_{s}(x).
\tag{4}
\end{align*}

The first equality unravels the definition of the mapping $\mathrm{ip}^{(1,X)}$; the second equality applies the Curry-Howard mapping at a $(1,0)$-identity path on a variable; the third equality unravels the definition of the mapping $\eta^{(0,X)}$ at a variable; finally, the last equality follows from the universal property of $\eta^{(1,0)\sharp}$ introduced in Proposition~\ref{PEmb}.

For the other equation, note that, by Propositions~\ref{PCHHom} and~\ref{PIpHom}, $\mathrm{CH}^{(1)}\circ\mathrm{ip}^{(1,0)\sharp}$ is a $\Sigma$-homomorphism. Moreover, the following chain of equalities holds
\allowdisplaybreaks
\begin{align*}
\mathrm{CH}^{(1)}\circ\mathrm{ip}^{(1,0)\sharp}\circ\eta^{(0,X)} &=
\mathrm{CH}^{(1)}\circ\mathrm{ip}^{(1,X)}
\tag{1}
\\&=
\eta^{(1,X)}.
\tag{2}
\end{align*}

The first equation follows from the universal property of $\mathrm{ip}^{(1,0)\sharp}$ introduced in Proposition~\ref{PIpHom}; the last equation follows from the equation proved before. 

Taking into account the unicity of the $\Sigma$-homomorphism $\eta^{(1,0)\sharp}$ from Proposition~\ref{PEmb}, we infer that $\mathrm{CH}^{(1)}\circ\mathrm{ip}^{(1,0)\sharp}=\eta^{(1,0)\sharp}$.
\end{proof}

Let us note that the Curry-Howard mapping is also an $S$-sorted mapping from the underlying $S$-sorted set of the  partial $\Sigma^{\boldsymbol{\mathcal{A}}}$-algebra $\mathbf{Pth}_{\boldsymbol{\mathcal{A}}}$ to the underlying $S$-sorted set of the free $\Sigma^{\boldsymbol{\mathcal{A}}}$-algebra $\mathbf{T}_{\Sigma^{\boldsymbol{\mathcal{A}}}}(X)$. The question of whether $\mathrm{CH}^{(1)}$ is a $\Sigma^{\boldsymbol{\mathcal{A}}}$-homomorphism immediately raises. In this regard, we next prove that this mapping preserves, at least, the constants from the $S$-sorted set $\mathcal{A}$.

\begin{figure}
\begin{center}
\begin{tikzpicture}
[ACliment/.style={-{To [angle'=45, length=5.75pt, width=4pt, round]}
}, scale=0.8]
\node[] (a) 		at 	(0,0) 	[] 	{$\mathcal{A}$};
\node[] (p) 		at 	(6,0) 	[] 	{$\mathrm{Pth}_{\boldsymbol{\mathcal{A}}}$};
\node[] (txqc) 	at 	(6,-3) 	[] 	{$\mathrm{T}_{\Sigma^{\boldsymbol{\mathcal{A}}}}
(X)$};
\draw[ACliment, bend right=10]   (a) 	to node [below left]	
{$\eta^{(1,\mathcal{A})}$} (txqc);
\draw[ACliment] (a) 	to node [above right]	
{$\mathrm{ech}^{(1,\mathcal{A})}$} (p);
\draw[ACliment]  (p) 	to node [right]	
{$\mathrm{CH}^{(1)}$} (txqc);
\end{tikzpicture}
\end{center}
\caption{Behaviour of $\mathrm{CH}^{(1)}$ relative to $\mathcal{A}$ at layer 1.}
\label{FCHA}
\end{figure}

\begin{proposition}\label{PCHA} For the Curry-Howard mapping we have that
$$
\mathrm{CH}^{(1)}\circ\mathrm{ech}^{(1,\mathcal{A})}=\eta^{(1,\mathcal{A})}.
$$
\end{proposition}

\begin{proof}
The reader is advised to consult the diagram presented in Figure~\ref{FCHA}.

Let $s$ be a sort in $S$ and $\mathfrak{p}$ a rewrite rule in $\mathcal{A}_{s}$, then the following chain of equalities holds
\begin{align*}
\mathrm{CH}^{(1)}_{s}\left(
\mathrm{ech}^{(1,\mathcal{A})}_{s}(\mathfrak{p})
\right)
&=
\mathrm{CH}^{(1)}_{s}\left(
\mathfrak{p}^{\mathbf{Pth}_{\boldsymbol{\mathcal{A}}}}
\right)
\tag{1}
\\&=
\mathfrak{p}^{\mathbf{T}_{\Sigma^{\boldsymbol{\mathcal{A}}}}(X)}
\tag{2}
\\&=
\eta^{(1,\mathcal{A})}_{s}(\mathfrak{p}).
\tag{3}
\end{align*}

The first equality recovers the interpretation of the constant symbol $\mathfrak{p}$ in the partial $\Sigma^{\boldsymbol{\mathcal{A}}}$-algebra $\mathbf{Pth}_{\boldsymbol{\mathcal{A}}}$ introduced in Proposition~\ref{PPthCatAlg}; the second equality follows by Definition~\ref{DCH}, since the Curry-Howard mapping assigns to an echelon path the unique rewrite rule appearing in it; finally, the last equality recovers the interpretation of the mapping $\eta^{(1,\mathcal{A})}$ at a rewrite rule $\mathfrak{p}$
as it is presented in Definition~\ref{DEta}.

This finishes the proof.
\end{proof}

However the preservation of  the constants from $\mathcal{A}$ is not enough to infer that the Curry-Howard mapping is a  $\Sigma^{\boldsymbol{\mathcal{A}}}$-homomorphism. This is due, essentially, to the fact that several terms can denote the same path. We next present a counterexample showing that $\mathrm{CH}^{(1)}$ is not a $\Sigma^{\boldsymbol{\mathcal{A}}}$-homomorphism.

\begin{restatable}{proposition}{PCHNotHomCat}
\label{PCHNotHomCat}
The Curry-Howard mapping $\mathrm{CH}^{(1)}$ is not a $\Sigma^{\boldsymbol{\mathcal{A}}}$-homomorphism from the partial $\Sigma^{\boldsymbol{\mathcal{A}}}$-algebra $\mathbf{Pth}_{\boldsymbol{\mathcal{A}}}$ to the free $\Sigma^{\boldsymbol{\mathcal{A}}}$-algebra $\mathbf{T}_{\Sigma^{\boldsymbol{\mathcal{A}}}}(X)$.
\end{restatable}
\begin{proof}
Let $s$ be a sort in $S$ and $P$ a term in $\mathrm{T}_{\Sigma}(X)_{s}$. Let us consider the $(1,0)$-identity path on $P$, i.e., $\mathrm{ip}^{(1,0)\sharp}_{s}(P)$. In virtue of Proposition~\ref{PCHId}, we have that
$$
\mathrm{CH}^{(1)}_{s}\left(
\mathrm{ip}^{(1,0)\sharp}_{s}
\left(
P
\right)\right)=\eta^{(1,0)\sharp}_{s}\left(P\right).
$$
However, the $(1,0)$-identity path on $P$ can be $0$-composed with itself and we have 
$$
\mathrm{ip}^{(1,0)\sharp}_{s}(P)
\circ_{s}^{0\mathbf{Pth}_{\boldsymbol{\mathcal{A}}}}
\mathrm{ip}^{(1,0)\sharp}_{s}(P)=\mathrm{ip}^{(1,0)\sharp}_{s}(P).
$$

Let us assume, towards a contradiction, that $\mathrm{CH}^{(1)}$ is a $\Sigma^{\boldsymbol{\mathcal{A}}}$-homomorphism. Then we have the following chain of equalities
\begin{align*}
\eta^{(1,0)\sharp}_{s}(P)
&=
\mathrm{CH}^{(1)}_{s}\left(
\mathrm{ip}^{(1,0)\sharp}_{s}\left(
P\right)
\right)
\tag{1}
\\&=
\mathrm{CH}^{(1)}_{s}\left(
\mathrm{ip}^{(1,0)\sharp}_{s}\left(P\right)
\circ_{s}^{0\mathbf{Pth}_{\boldsymbol{\mathcal{A}}}}
\mathrm{ip}^{(1,0)\sharp}_{s}\left(P\right)
\right)
\tag{2}
\\&=
\mathrm{CH}^{(1)}_{s}\left(
\mathrm{ip}^{(1,0)\sharp}_{s}(P)
\right)
\circ_{s}^{0\mathbf{T}_{\Sigma^{\boldsymbol{\mathcal{A}}}}(X)}
\mathrm{CH}^{(1)}_{s}\left(
\mathrm{ip}^{(1,0)\sharp}_{s}(P)
\right)
\tag{3}
\\&=
\eta^{(1,0)\sharp}_{s}(P)
\circ_{s}^{0\mathbf{T}_{\Sigma^{\boldsymbol{\mathcal{A}}}}(X)}
\eta^{(1,0)\sharp}_{s}(P).
\tag{4}
\end{align*}

The first equality follows from Proposition~\ref{PCHId}; the second equality follows from the fact that $\mathrm{ip}^{(1,0)\sharp}_{s}(P)$ is idempotent for the $0$-composition; the third equality follows since we are assuming that $\mathrm{CH}^{(1)}$ is a $\Sigma^{\boldsymbol{\mathcal{A}}}$-homomorphism; finally the last equality follows from Proposition~\ref{PCHId}.

However, $\eta^{(1,0)\sharp}_{s}(P)\circ_{s}^{0\mathbf{T}_{\Sigma^{\boldsymbol{\mathcal{A}}}}(X)}\eta^{(1,0)\sharp}_{s}(P)$ and $\eta^{(1,0)\sharp}_{s}(P)$ are different terms in $\mathrm{T}_{\Sigma^{\boldsymbol{\mathcal{A}}}}(X)_{s}$.
\end{proof}

Nevertheless, there are some specific configurations on  $0$-compositions of paths for which the  Curry-Howard mapping behaves like a  $\Sigma^{\boldsymbol{\mathcal{A}}}$-homomorphism. As we will make use of these configurations later, we specify them in the following proposition.

\begin{proposition}\label{PCHHomCat}
Let $s$ be a sort in $S$ and let $\mathfrak{P}$ be a path in $\mathrm{Pth}_{\boldsymbol{\mathcal{A}},s}$ of length strictly greater than one containing at least one echelon. Let $i\in\bb{\mathfrak{P}}$ be the first index for which $\mathfrak{P}^{i,i}$ is an echelon. Then
\begin{itemize}
\item[(i)] if $i=0$, that is, if the first echelon appears in the first step, then
$$
\mathrm{CH}^{(1)}_{s}
\left(
\mathfrak{P}^{1,\bb{\mathfrak{P}}-1}
\circ_{s}^{0\mathbf{Pth}_{\boldsymbol{\mathcal{A}}}}
\mathfrak{P}^{0,0}
\right)
=
\mathrm{CH}^{(1)}_{s}
\left(
\mathfrak{P}^{1,\bb{\mathfrak{P}}-1}
\right)
\circ_{s}^{0\mathbf{T}_{\Sigma^{\boldsymbol{\mathcal{A}}}}(X)}
\mathrm{CH}^{(1)}_{s}
\left(
\mathfrak{P}^{0,0}
\right);
$$
\item[(ii)] if $i=\bb{\mathfrak{P}}-1$, that is, if the first echelon appears in the last step, then
\begin{multline*}
\mathrm{CH}^{(1)}_{s}
\left(
\mathfrak{P}^{\bb{\mathfrak{P}}-1,\bb{\mathfrak{P}}-1}
\circ_{s}^{0\mathbf{Pth}_{\boldsymbol{\mathcal{A}}}}
\mathfrak{P}^{0,\bb{\mathfrak{P}}-2}
\right)
\\=
\mathrm{CH}^{(1)}_{s}
\left(
\mathfrak{P}^{\bb{\mathfrak{P}}-1,\bb{\mathfrak{P}}-1}
\right)
\circ_{s}^{0\mathbf{T}_{\Sigma^{\boldsymbol{\mathcal{A}}}}(X)}
\mathrm{CH}^{(1)}_{s}
\left(
\mathfrak{P}^{0,\bb{\mathfrak{P}}-2}
\right).
\end{multline*}
\end{itemize}
\end{proposition}
\begin{proof}
These cases were already considered in Definition~\ref{DCH}.
\end{proof}

In the following example we compute the value of the Curry-Howard mapping at the path introduced in Example~\ref{ERun}.

\begin{example}[Continuation of Example~\ref{ERun}]\label{ERunII}
The categorial signature determined by $\mathcal{A}$ on $\Sigma$, i.e., $\Sigma^{\boldsymbol{\mathcal{A}}}$, is the $2$-sorted signature defined, for every $(\mathbf{s},s)\in 2^{\star}\times 2$, as
$$
\Sigma^{\boldsymbol{\mathcal{A}}}_{\mathbf{s},s}=
\left\{
\begin{array}{llll}
\{\alpha\}\amalg \{\mathfrak{p}_{k}\mid k\in 8\},&\mbox{if }\mathbf{s}=\lambda \mbox{ and } s=0;\\
 \{\mathfrak{q}_{k}\mid k\in 2\},&\mbox{if }\mathbf{s}=\lambda \mbox{ and } s=1;\\
\{\mathrm{sc}^{0}_{0}, \mathrm{tg}^{0}_{0}\},&\mbox{if }\mathbf{s}=(0) \mbox{ and } s=0;\\
\{\circ^{0}_{0}\},&\mbox{if }\mathbf{s}=(0,0) \mbox{ and } s=0;\\
\{\beta\},&\mbox{if }\mathbf{s}=(0,0,0) \mbox{ and }=0;\\
\{\mathrm{sc}^{0}_{1}, \mathrm{tg}^{0}_{1}\},&\mbox{if }\mathbf{s}=(1) \mbox{ and } s=1;\\
\{\gamma\},&\mbox{if }\mathbf{s}=(1,0) \mbox{ and } s=0;\\
\{\delta\},&\mbox{if }\mathbf{s}=(1,1) \mbox{ and } s=0;\\
\{\circ^{0}_{1}\},&\mbox{if }\mathbf{s}=(1,1) \mbox{ and } s=1;\\
\varnothing&\mbox{otherwise.}
\end{array}\right.
$$

We can now compute the value of the Curry-Howard mapping at $\mathfrak{P}$. To simplify the presentation we have omitted  the superscript $\mathbf{T}_{\Sigma^{\boldsymbol{\mathcal{A}}}}(X)$.

First, we note that $\mathfrak{P}$ contains echelons and that, actually, the first echelon appears at the fourth step. Therefore, since $\mathfrak{P}^{4,4}$ is the echelon associated to $\mathfrak{p}_{6}$, $\mathrm{CH}^{(1)}_{0}(\mathfrak{P}^{4,4})=\mathfrak{p}_{6}$. Thus 
\allowdisplaybreaks
\begin{align*}
\mathrm{CH}^{(1)}_{0}\left(
\mathfrak{P}\right)&=\mathrm{CH}^{(1)}_{0}\left(
\mathfrak{P}^{4,13}\right)\circ^{0}_{0}\mathrm{CH}^{(1)}_{0}
\left(
\mathfrak{P}^{0,3}
\right)\\
&=
\left(\mathrm{CH}^{(1)}_{0}\left(\mathfrak{P}^{5,13}\right)\circ^{0}_{0}\mathfrak{p}_{6}\right)\circ^{0}_{0}\mathrm{CH}^{(1)}_{0}\left(\mathfrak{P}^{0,3}\right).
\end{align*}

But, by Example~\ref{ERun}, $\mathfrak{P}^{0,3}$ is an echelonless path for which the extraction algorithm from Lemma~\ref{LPthExtract} retrieves the family of paths $(\mathfrak{Q}_{j})_{j\in 3}$, which is composed of echelons or compositions of echelons. Hence, the value of the Curry-Howard mapping at $\mathfrak{P}^{0,3}$ is
\begin{align*}
\mathrm{CH}^{(1)}_{0}\left(
\mathfrak{P}^{0,3}
\right)&=\beta\left(\left(\mathrm{CH}^{(1)}_{0}\left(\mathfrak{Q}_{j}\right)\right)_{j\in 3}\right)=\beta\left(\mathfrak{p}_{4},\mathfrak{p}_{1}\circ^{0}_{0}\mathfrak{p}_{2},\mathfrak{p}_{1}\right).
\end{align*}

Moreover, the value of the Curry-Howard mapping at the subpath $\mathfrak{P}^{5,13}$ of $\mathfrak{P}$ is
\begin{flushleft}
$\mathrm{CH}^{(1)}_{0}\left(\mathfrak{P}^{5,13}\right)$
\allowdisplaybreaks
\begin{align*}
\qquad&=
\mathrm{CH}^{(1)}_{0}\left(
\mathfrak{P}^{7,13}
\right)
\circ_{0}^{0}
\mathrm{CH}^{(1)}_{0}\left(
\mathfrak{P}^{5,6}
\right)
\\&=
\left(
\mathrm{CH}^{(1)}_{0}\left(
\mathfrak{P}^{8,13}
\right)
\circ^{0}_{0}
\mathfrak{p}_{5}
\right)
\circ_{0}^{0}
\delta\left(x,\mathfrak{q}_{0}\circ_{1}^{0}\mathfrak{q}_{1}\right)
\\&=
\left(
\left(
\mathrm{CH}^{(1)}_{0}\left(
\mathfrak{P}^{12,13}
\right)
\circ^{0}_{0}
\mathrm{CH}^{(1)}_{0}\left(
\mathfrak{P}^{8,11}
\right)
\right)
\circ^{0}_{0}
\mathfrak{p}_{5}
\right)
\circ_{0}^{0}
\delta\left(x,\mathfrak{q}_{0}\circ_{1}^{0}\mathfrak{q}_{1}\right)
\\&=
\left(
\left(
\left(
\mathrm{CH}^{(1)}_{0}\left(
\mathfrak{P}^{13,13}
\right)
\circ^{0}_{0}
\mathrm{CH}^{(1)}_{0}\left(
\mathfrak{P}^{12,12}
\right)
\right)
\circ^{0}_{0}
\beta\left(
\mathfrak{p}_{4},
\delta\left(\mathfrak{q}_{0},\mathfrak{q}_{1}\right)
\circ^{0}_{0}
\mathfrak{p}_{3},
a
\right)
\right)
\circ^{0}_{0}
\mathfrak{p}_{5}
\right)
\\&
\qquad\qquad\qquad\qquad\qquad\qquad\qquad
\qquad\qquad\qquad\qquad\qquad\qquad
\circ_{0}^{0}
\delta\left(x,\mathfrak{q}_{0}\circ_{1}^{0}\mathfrak{q}_{1}\right)
\\&=
\left(
\left(
\left(
\mathfrak{p}_{0}
\circ^{0}_{0}
\mathfrak{p}_{7}
\right)
\circ^{0}_{0}
\beta\left(
\mathfrak{p}_{4},
\delta\left(\mathfrak{q}_{0},\mathfrak{q}_{1}\right)
\circ^{0}_{0}
\mathfrak{p}_{3},
a
\right)
\right)
\circ^{0}_{0}
\mathfrak{p}_{5}
\right)
\circ_{0}^{0}
\delta\left(x,\mathfrak{q}_{0}\circ_{1}^{0}\mathfrak{q}_{1}\right).
\end{align*}
\end{flushleft}

Therefore, the value of the Curry-Howard mapping at $\mathfrak{P}$ is 
\begin{multline*}
\mathrm{CH}^{(1)}_{0}\left(\mathfrak{P}\right)=
\Big(
\left(
\left(
\left(
\mathfrak{p}_{0}
\circ^{0}_{0}
\mathfrak{p}_{7}
\right)
\circ^{0}_{0}
\beta\left(
\mathfrak{p}_{4},
\delta\left(\mathfrak{q}_{0},\mathfrak{q}_{1}\right)
\circ^{0}_{0}
\mathfrak{p}_{3},
a
\right)
\right)
\circ^{0}_{0}
\mathfrak{p}_{5}
\right)
\circ_{0}^{0}
\\
\delta\left(x,\mathfrak{q}_{0}\circ_{1}^{0}\mathfrak{q}_{1}\right)
\circ^{0}_{0}\mathfrak{p}_{6}
\Big)
\circ^{0}_{0}
\beta\left(\mathfrak{p}_{4},\mathfrak{p}_{1}\circ^{0}_{0}\mathfrak{p}_{2},\mathfrak{p}_{1}\right).
\end{multline*}

This is the computation of the value of the Curry-Howard mapping at $\mathfrak{P}$.
\end{example}

We conclude this section by proving that $\coprod\mathrm{CH}^{(1)}$ is a monotone mapping from $(\coprod\mathrm{Pth}_{\boldsymbol{\mathcal{A}}}, \leq_{\mathbf{Pth}_{\boldsymbol{\mathcal{A}}}})$ to 
$(\coprod\mathrm{T}_{\Sigma^{\boldsymbol{\mathcal{A}}}}(X), \leq_{\mathbf{T}_{\Sigma^{\boldsymbol{\mathcal{A}}}}(X)})$. Ultimately, this will be crucial as the more complex the paths, the more complex the associated terms will be.

\begin{restatable}{proposition}{PCHMono}
\label{PCHMono} 
The mapping $\coprod\mathrm{CH}^{(1)}$ from
$\coprod\mathrm{Pth}_{\boldsymbol{\mathcal{A}}}$ to  $\coprod\mathrm{T}_{\Sigma^{\boldsymbol{\mathcal{A}}}}(X)$ is order-preserving
$$
\textstyle
\coprod\mathrm{CH}^{(1)}\colon
\left(
\coprod\mathrm{Pth}_{\boldsymbol{\mathcal{A}}}, \leq_{\mathbf{Pth}_{\boldsymbol{\mathcal{A}}}}
\right)
\mor
\left(
\coprod\mathrm{T}_{\Sigma^{\boldsymbol{\mathcal{A}}}}(X), \leq_{\mathbf{T}_{\Sigma^{\boldsymbol{\mathcal{A}}}}(X)}
\right),
$$
that is, given pairs $(\mathfrak{Q},t),(\mathfrak{P},s)$ in $\coprod\mathrm{Pth}_{\boldsymbol{\mathcal{A}}}$, if $(\mathfrak{Q},t)\leq_{\mathbf{Pth}_{\boldsymbol{\mathcal{A}}}} (\mathfrak{P},s)$, then $\mathrm{CH}^{(1)}_{t}(\mathfrak{Q})$ is a subterm of type $t$ of the term $\mathrm{CH}^{(1)}_{s}(\mathfrak{P})$, i.e., 
$$
\mathrm{CH}^{(1)}_{t}\left(
\mathfrak{Q}
\right)\in\mathrm{Subt}_{\Sigma^{\boldsymbol{\mathcal{A}}}}\left(
\mathrm{CH}^{(1)}_{s}\left(
\mathfrak{P}
\right)
\right)_{t}.
$$
\end{restatable}

\begin{proof}
Let us recall from Remark~\ref{ROrd} that $(\mathfrak{Q},t)\leq_{\mathbf{Pth}_{\boldsymbol{\mathcal{A}}}}(\mathfrak{P},s)$ if and only if $s=t$ and $\mathfrak{Q}=\mathfrak{P}$ or there exists a natural number $m\in\mathbb{N}-\{0\}$, a word $\mathbf{w}\in S^{\star}$ of length $\bb{\mathbf{w}}=m+1$, and a family of paths $(\mathfrak{R}_{k})_{k\in\bb{\mathbf{w}}}$ in $\mathrm{Pth}_{\boldsymbol{\mathcal{A}},\mathbf{w}}$ such that $w_{0}=t$, $\mathfrak{R}_{0}=\mathfrak{Q}$, $w_{m}=s$, $\mathfrak{R}_{m}=\mathfrak{P}$ and, for every $k\in m$, $(\mathfrak{R}_{k}, w_{k})\prec_{\mathbf{Pth}_{\boldsymbol{\mathcal{A}}}} (\mathfrak{R}_{k+1}, w_{k+1})$.

The proposition holds trivially for $s=t$ and $\mathfrak{Q}=\mathfrak{P}$. Therefore, it remains to prove the other case.  We will prove it by induction on $m\in\mathbb{N}-\{0\}$.

\textsf{Base step of the induction.}

For $m=1$ we have that $(\mathfrak{Q},t)\prec_{\mathbf{Pth}_{\boldsymbol{\mathcal{A}}}}(\mathfrak{P},s)$, hence following Definition~\ref{DOrd} we are in one of the following situations
\begin{enumerate}
\item $\mathfrak{P}$ is a $(1,0)$-identity path; or
\item $\mathfrak{P}$ is a path of length strictly greater than one containing at least one echelon; or
\item $\mathfrak{P}$ is an echelonless path.
\end{enumerate}

If~(1), then $\mathfrak{P}$ and $\mathfrak{Q}$ are $(1,0)$-identity paths of the form
\begin{align*}
\mathfrak{P}&=
\mathrm{ip}^{(1,0)\sharp}_{s}(P),
&
\mathfrak{Q}&=
\mathrm{ip}^{(1,0)\sharp}_{t}(Q),
\end{align*}
for some terms $P\in\mathrm{T}_{\Sigma}(X)_{s}$ and $Q\in\mathrm{T}_{\Sigma}(X)_{t}$ and the following inequality holds
$$
(Q,t)
<_{\mathbf{T}_{\Sigma}(X)}
(P,s)
$$
where $\leq_{\mathbf{T}_{\Sigma}(X)}$ is the Artinian partial order on $\coprod\mathrm{T}_{\Sigma}(X)$ introduced in Remark~\ref{RTermOrd}. Therefore, $Q$ is a subterm of $P$ of type $t$, that is $Q\in\mathrm{Subt}_{\Sigma}(P)_{t}$.

Note that, from Definition~\ref{DCH}, the following equalities hold
\begin{align*}
\mathrm{CH}^{(1)}_{t}\left(
\mathfrak{Q}
\right)&=
\eta^{(1,0)\sharp}_{t}(Q),
&
\mathrm{CH}^{(1)}_{s}\left(
\mathfrak{P}
\right)&=
\eta^{(1,0)\sharp}_{t}(P).
\end{align*}

Since $\coprod\eta^{(1,0)\sharp}$ is an order embedding of $(\coprod\mathrm{T}_{\Sigma}(X),\leq_{\mathbf{T}_{\Sigma}(X)})$ into the partially ordered set $(\coprod\mathrm{T}_{\Sigma^{\boldsymbol{\mathcal{A}}}}(X),\leq_{\mathbf{T}_{\Sigma^{\boldsymbol{\mathcal{A}}}}(X)})$, we have that 
$$
\left(
\mathrm{CH}^{(1)}_{t}\left(
\mathfrak{Q}
\right),
t
\right)
\leq_{\mathbf{T}_{\Sigma^{\boldsymbol{\mathcal{A}}}}(X)}
\left(
\mathrm{CH}^{(1)}_{s}\left(
\mathfrak{P}
\right),
s
\right).
$$

Case~(1) follows.

If~(2), i.e., if $\mathfrak{P}$ is a path of length strictly greater than one containing its first echelon at position $i\in\bb{\mathfrak{P}}$ then, depending on this position, we have the following subcases.

If $i=0$, then according to Definition~\ref{DOrd}, $\mathfrak{Q}$ is equal to $\mathfrak{P}^{0,0}$ or $\mathfrak{P}^{1,\bb{\mathfrak{P}}-1}$ and, according to Definition~\ref{DCH}, we have that
$$
\mathrm{CH}^{(1)}_{s}\left(
\mathfrak{P}
\right)=
\mathrm{CH}^{(1)}_{s}
\left(
\mathfrak{P}^{1,\bb{\mathfrak{P}}-1}
\right)
\circ_{s}^{0\mathbf{T}_{\Sigma^{\boldsymbol{\mathcal{A}}}}(X)}
\mathrm{CH}^{(1)}_{s}\left(
\mathfrak{P}^{0,0}
\right).
$$
In any case, we have that $(\mathrm{CH}^{(1)}_{s}(\mathfrak{Q}),s)
\leq_{\mathbf{T}_{\Sigma^{\boldsymbol{\mathcal{A}}}}(X)}
(\mathrm{CH}^{(1)}_{s}(\mathfrak{P}),s)
$.

If $i> 0$, then according to Definition~\ref{DOrd}, $\mathfrak{Q}$ is equal to $\mathfrak{P}^{0,i-1}$ or $\mathfrak{P}^{i,\bb{\mathfrak{P}}-1}$ and, according to Definition~\ref{DCH}, we have that
$$
\mathrm{CH}^{(1)}_{s}\left(
\mathfrak{P}
\right)=
\mathrm{CH}^{(1)}_{s}
\left(
\mathfrak{P}^{i,\bb{\mathfrak{P}}-1}
\right)
\circ_{s}^{0\mathbf{T}_{\Sigma^{\boldsymbol{\mathcal{A}}}}(X)}
\mathrm{CH}^{(1)}_{s}\left(
\mathfrak{P}^{0,i-1}
\right).
$$
In any case, we have that $(\mathrm{CH}^{(1)}_{s}(\mathfrak{Q}),s)
\leq_{\mathbf{T}_{\Sigma^{\boldsymbol{\mathcal{A}}}}(X)}
(\mathrm{CH}^{(1)}_{s}(\mathfrak{P}),s)
$.

Case~(2) follows.

If~(3), i.e., if $\mathfrak{P}$ is an echelonless path then, according to Lemma~\ref{LPthHeadCt},  there exists a unique word $\mathbf{s}\in S^{\star}-\{\lambda\}$ and a unique operation symbol $\sigma\in \Sigma_{\mathbf{s},s}$ associated to $\mathfrak{P}$. Let $(\mathfrak{P}_{j})_{j\in\bb{\mathbf{s}}}$ be the family of paths in $\mathrm{Pth}_{\boldsymbol{\mathcal{A}},\mathbf{s}}$ which, in virtue of Lemma~\ref{LPthExtract}, we can extract from $\mathfrak{P}$.  Then, according to Definition~\ref{DOrd}, $\mathfrak{Q}$ is equal to $\mathfrak{P}_{j}$, for some $j\in\bb{\mathbf{s}}$, and, according to Definition~\ref{DCH}, we have that 
$$
\mathrm{CH}^{(1)}_{s}
\left(
\mathfrak{P}
\right)=
\sigma^{\mathbf{T}_{\Sigma^{\boldsymbol{\mathcal{A}}}}(X)}
\left(\left(
\mathrm{CH}^{(1)}_{s_{j}}\left(
\mathfrak{P}_{j}
\right)\right)_{j\in\bb{\mathbf{s}}}
\right).
$$

In any case, we have that 
$
(\mathrm{CH}^{(1)}_{s_{j}}(\mathfrak{Q}),s_{j})
\leq_{\mathbf{T}_{\Sigma^{\boldsymbol{\mathcal{A}}}}(X)}
(\mathrm{CH}^{(1)}_{s}(\mathfrak{P}),s).
$

This concludes Case~(3).

This completes the base step.

\textsf{Inductive step of the induction.}

Assume the statement holds for sequences of length up to $m$, i.e., for every pair of sorts $t,s\in S$, if $(\mathfrak{Q},t), (\mathfrak{P},s)$ are pairs in $\coprod
\mathrm{Pth}_{\boldsymbol{\mathcal{A}}}$ such that  there exists a word $\mathbf{w}\in S^{\star}$ of length $\bb{\mathbf{w}}=m+1$ and a family of paths $(\mathfrak{R}_{k})_{k\in\bb{\mathbf{w}}}$ in $\mathrm{Pth}_{\boldsymbol{\mathcal{A}},\mathbf{w}}$  such that $w_{0}=t$, $\mathfrak{R}_{0}=\mathfrak{Q}$, $w_{m}=s$, $\mathfrak{R}_{m}=\mathfrak{P}$ and 
for every $k\in m$, $(\mathfrak{R}_{k}, w_{k})\prec_{\mathbf{Pth}_{\boldsymbol{\mathcal{A}}}} (\mathfrak{R}_{k+1}, w_{k+1})$ then 
$
\mathrm{CH}^{(1)}_{t}(\mathfrak{Q})\in\mathrm{Subt}_{\Sigma^{\boldsymbol{\mathcal{A}}}}(\mathrm{CH}^{(1)}_{s}(\mathfrak{P}))_{t}.
$

Now we prove it for sequences of length $m+1$. 

Let $t,s$ be sorts in $S$, if $(\mathfrak{Q},t), (\mathfrak{P},s)$ are pairs in $\coprod
\mathrm{Pth}_{\boldsymbol{\mathcal{A}}}$ such that  there exists a word $\mathbf{w}\in S^{\star}$ of length $\bb{\mathbf{w}}=m+2$ and a family of paths $(\mathfrak{R}_{k})_{k\in\bb{\mathbf{w}}}$ in $\mathrm{Pth}_{\boldsymbol{\mathcal{A}},\mathbf{w}}$  such that $w_{0}=t$, $\mathfrak{R}_{0}=\mathfrak{Q}$, $w_{m+1}=s$, $\mathfrak{R}_{m+1}=\mathfrak{P}$ and 
for every $k\in m+1$, $(\mathfrak{R}_{k}, w_{k})\prec_{\mathbf{Pth}_{\boldsymbol{\mathcal{A}}}} (\mathfrak{R}_{k+1}, w_{k+1})$.

Consider the word $\mathbf{w}^{0,m}$ of length $\bb{\mathbf{w}^{0,m}}=m+1$ and the family of paths $(\mathfrak{R}_{k})_{k\in\bb{\mathbf{w}^{0,m}}}$ in $\mathrm{Pth}_{\boldsymbol{\mathcal{A}},\mathbf{w}^{0,m}}$. This is a sequence of length $m$ instantiating that 
$
(\mathfrak{Q},t)
\leq_{\mathbf{Pth}_{\boldsymbol{\mathcal{A}}}}
(\mathfrak{R}_{m}, w_{m})$. By the inductive hypothesis, we have that 
$\mathrm{CH}^{(1)}_{t}(\mathfrak{Q})
\in\mathrm{Subt}_{\Sigma^{\boldsymbol{\mathcal{A}}}}(\mathrm{CH}^{(1)}_{w_{m}}(\mathfrak{R}_{m}))_{t}.
$

Now, consider the sequence of paths $(\mathfrak{R}_{m},\mathfrak{P})$ it is a one-step sequence of paths in $\mathrm{Pth}_{\boldsymbol{\mathcal{A}},\mathbf{w}^{m,m+1}}$ satisfying that $(\mathfrak{R}_{m},w_{m})\prec_{\mathbf{Pth}_{\boldsymbol{\mathcal{A}}}} (\mathfrak{P}, s)$. By the base case, we have that 
$
\mathrm{CH}^{(1)}_{w_{m}}(\mathfrak{R}_{m})
\in\mathrm{Subt}_{\Sigma^{\boldsymbol{\mathcal{A}}}}(\mathrm{CH}^{(1)}_{s}(\mathfrak{P}))_{w_{m}}.
$

Hence, by the transitivity of the partial order $\leq_{\mathbf{T}_{\Sigma^{\boldsymbol{\mathcal{A}}}}(X)}$, we have that 
$$\mathrm{CH}^{(1)}_{t}
\left(\mathfrak{Q}
\right)
\in\mathrm{Subt}_{\Sigma^{\boldsymbol{\mathcal{A}}}}
\left(\mathrm{CH}^{(1)}_{s}
\left(\mathfrak{P}
\right)\right)_{t}.
$$

This finishes the proof.
\end{proof}

\chapter{
\texorpdfstring
{On the kernel of the Curry-Howard mapping}
{The kernel of the Curry-Howard mapping}
}\label{S1E}

In this chapter we revisit the Curry-Howard mapping. Specifically, our objective will be to prove that 
$\mathrm{Ker}(\mathrm{CH}^{(1)})$, the kernel of $\mathrm{CH}^{(1)}$, is a closed $\Sigma^{\boldsymbol{\mathcal{A}}}$-congruence on $\mathbf{Pth}_{\boldsymbol{\mathcal{A}}}$. To achieve such an objective it will be necessary to prove several auxiliary results. So, we first prove that the Curry-Howard mapping is able to completely characterise $(1,0)$-identity paths. In fact, we prove that $(1,0)$-identity paths in 
$\mathrm{Ker}(\mathrm{CH}^{(1)})$ must be necessarily equal. The Curry-Howard mapping is also capable of characterising echelons. We prove that paths in $\mathrm{Ker}(\mathrm{CH}^{(1)})$ necessarily have the same length, the same $(0,1)$-source and $(0,1)$-target. This, in fact, helps to demonstrate that paths of length one are completely characterised by the Curry-Howard mapping. Indeed, two paths of length one in 
$\mathrm{Ker}(\mathrm{CH}^{(1)})$ must be necessarily equal. After having done all that we prove that 
$\mathrm{Ker}(\mathrm{CH}^{(1)})$ is a closed $\Sigma^{\boldsymbol{\mathcal{A}}}$-congruence on 
$\mathbf{Pth}_{\boldsymbol{\mathcal{A}}}$.


To begin with, we identify $\mathcal{A}$ with a subset of terms in $\mathrm{T}_{\Sigma^{\boldsymbol{\mathcal{A}}}}(X)$ by means of its image under the $S$-sorted mapping $\eta^{(1,\mathcal{A})}$ presented in Definition~\ref{DEta}, that is
$$
\eta^{(1,\mathcal{A})}
\left[\mathcal{A}\right]=
\left(\left\lbrace
\mathfrak{p}^{\mathbf{T}_{\Sigma^{\boldsymbol{\mathcal{A}}}}(X)}\in 
\mathrm{T}_{\Sigma^{\boldsymbol{\mathcal{A}}}}(X)_{s}
\,\middle|
\,
\mathfrak{p}\in\mathcal{A}_{s}
\right\rbrace\right)_{s\in S}.
$$

We next prove that the value of the Curry-Howard mapping at a non-$(1,0)$-identity path always has a subterm which contains a rewrite rule (of some sort) in $\mathcal{A}$. 

\begin{proposition}\label{PCHRewId}
Let $s$ be a sort in $S$ and $\mathfrak{P}$ a non-$(1,0)$-identity path in $\mathrm{Pth}_{\boldsymbol{\mathcal{A}},s}$. Then 
$$
\eta^{(1,\mathcal{A})}
\left[
\mathcal{A}
\right]\cap\mathrm{Subt}_{\Sigma^{\boldsymbol{\mathcal{A}}}}
\left(
\mathrm{CH}^{(1)}_{s}\left(
\mathfrak{P}
\right)\right)
\neq\varnothing^{S}.
$$
\end{proposition}
\begin{proof}
We prove the statement by Artinian induction on $(\coprod\mathrm{Pth}_{\boldsymbol{\mathcal{A}}},\leq_{\mathbf{Pth}_{\boldsymbol{\mathcal{A}}}})$.

\textsf{Base step of the Artinian induction}.

Let $(\mathfrak{P},s)$ be a minimal element of $(\coprod\mathrm{Pth}_{\boldsymbol{\mathcal{A}}}, 
\leq_{\mathbf{Pth}_{\boldsymbol{\mathcal{A}}}})$. Then, by Proposition~\ref{PMinimal}, the path $\mathfrak{P}$ is either~(1) a $(1,0)$-identity path on a minimal term or~(2) an echelon. Since we are assuming that $\mathfrak{P}$ is a non-$(1,0)$-identity path, we have that $\mathfrak{P}$ is an echelon.  Hence 
$\mathrm{CH}^{(1)}_{s}(\mathfrak{P})$ is $\mathfrak{p}^{\mathbf{T}_{\Sigma^{\boldsymbol{\mathcal{A}}}}(X)}$, for $\mathfrak{p}\in\mathcal{A}_{s}$, the rewrite rule defining the echelon. 
Therefore 
$$\mathfrak{p}^{\mathbf{T}_{\Sigma^{\boldsymbol{\mathcal{A}}}}(X)}
\in \eta^{(1,\mathcal{A})}
\left[
\mathcal{A}
\right]_{s}\cap\mathrm{Subt}_{\Sigma^{\boldsymbol{\mathcal{A}}}}
\left(
\mathrm{CH}^{(1)}_{s}\left(
\mathfrak{P}
\right)\right)_{s}.$$

\textsf{Inductive step of the Artinian induction}.

Let $(\mathfrak{P},s)$ be a non-minimal element of $(\coprod\mathrm{Pth}_{\boldsymbol{\mathcal{A}}}, \leq_{\mathbf{Pth}_{\boldsymbol{\mathcal{A}}}})$ such that $\mathfrak{P}$ is a non-$(1,0)$-identity path. Let us suppose that, for every sort $t\in S$ and every path $\mathfrak{Q}$ in $\mathrm{Pth}_{\boldsymbol{\mathcal{A}},t}$,  if $(\mathfrak{Q},t)<_{\mathbf{Pth}_{\boldsymbol{\mathcal{A}}}}$-precedes $(\mathfrak{P},s)$, then the statement holds for $\mathfrak{Q}$, i.e., if $\mathfrak{Q}$ a non-$(1,0)$-identity path in $\mathrm{Pth}_{\boldsymbol{\mathcal{A}},t}$, then 
$$
\eta^{(1,\mathcal{A})}\left[
\mathcal{A}
\right]\cap\mathrm{Subt}_{\Sigma^{\boldsymbol{\mathcal{A}}}}
\left(
\mathrm{CH}^{(1)}_{t}
\left(
\mathfrak{Q}
\right)\right)\neq\varnothing^{S}.
$$
Since $(\mathfrak{P},s)$ is a non-minimal element of $(\coprod\mathrm{Pth}_{\boldsymbol{\mathcal{A}}}, \leq_{\mathbf{Pth}_{\boldsymbol{\mathcal{A}}}})$, and taking into account that $\mathfrak{P}$ is a non-$(1,0)$-identity path, we have, by Lemma~\ref{LOrdI}, that $\mathfrak{P}$ is either~(1) a path of length strictly greater than one containing at least one echelon, or~(2) an echelonless path.

If~(1), then let $i\in \bb{\mathfrak{P}}$ be the first index for which the one-step subpath $\mathfrak{P}^{i,i}$ is an echelon.  In this case, according to the nature of $i$, we have that, if $i=0$, then $\mathfrak{P}^{0,0}$ is an echelon and, if $i>0$, then $\mathfrak{P}^{i,\bb{\mathfrak{P}}-1}$ contains an echelon.

Note that in any case we have found a subpath $\mathfrak{Q}$ of $\mathfrak{P}$ which is a non-$(1,0)$-identity path and such that $(\mathfrak{Q},s)$
$\prec_{\mathbf{Pth}_{\boldsymbol{\mathcal{A}}}}$-precedes $(\mathfrak{P},s)$.
Then, by induction, we have that
$
\eta^{(1,\mathcal{A})}[\mathcal{A}]\cap\mathrm{Subt}_{\Sigma^{\boldsymbol{\mathcal{A}}}}(\mathrm{CH}^{(1)}_{s}(\mathfrak{Q}))\neq\varnothing^{S}.
$ By Proposition~\ref{PCHMono}, we also have that $
\mathrm{CH}^{(1)}_{s}(\mathfrak{Q})\in\mathrm{Subt}_{\Sigma^{\boldsymbol{\mathcal{A}}}}(
\mathrm{CH}^{(1)}_{s}(\mathfrak{P})
)_{s}
$, hence 
$
\eta^{(1,\mathcal{A})}[\mathcal{A}]\cap\mathrm{Subt}_{\Sigma^{\boldsymbol{\mathcal{A}}}}(\mathrm{CH}^{(1)}_{s}(\mathfrak{P}))\neq\varnothing^{S}.
$

If~(2), i.e., if $\mathfrak{P}$ is an echelonless path then, in virtue of Lemma~\ref{LPthHeadCt},  there exists a unique word $\mathbf{s}\in S^{\star}-\{\lambda\}$ and a unique operation symbol $\sigma\in \Sigma_{\mathbf{s},s}$ associated to $\mathfrak{P}$. Let $(\mathfrak{P}_{j})_{j\in\bb{\mathbf{s}}}$ be the family of paths in $\mathrm{Pth}_{\boldsymbol{\mathcal{A}},\mathbf{s}}$ which, in virtue of Lemma~\ref{LPthExtract}, we can extract from $\mathfrak{P}$. Since $\mathfrak{P}$ is not an $(1,0)$-identity path, there exists at least one $j\in\bb{\mathbf{s}}$ such that $\mathfrak{P}_{j}$ is a non-$(1,0)$-identity path and, in addition, $(\mathfrak{P}_{j},s_{j})$ $\prec_{\mathbf{Pth}_{\boldsymbol{\mathcal{A}}}}$-precedes $(\mathfrak{P},s)$.
Then, by induction, we have that $\eta^{(1,\mathcal{A})}[\mathcal{A}]\cap\mathrm{Subt}_{\Sigma^{\boldsymbol{\mathcal{A}}}}(\mathrm{CH}^{(1)}_{s_{j}}(\mathfrak{P}_{j}))\neq\varnothing^{S}$.  By Proposition~\ref{PCHMono}, we also have that $
\mathrm{CH}^{(1)}_{s_{j}}(\mathfrak{P}_{j})\in\mathrm{Subt}_{\Sigma^{\boldsymbol{\mathcal{A}}}}(
\mathrm{CH}^{(1)}_{s}(\mathfrak{P})
)_{s_{j}}
$, hence
$
\eta^{(1,\mathcal{A})}[\mathcal{A}]\cap\mathrm{Subt}_{\Sigma^{\boldsymbol{\mathcal{A}}}}(\mathrm{CH}^{(1)}_{s}(\mathfrak{P}))\neq\varnothing^{S}.
$

This finishes the proof.
\end{proof}

In the following corollary we characterize the $(1,0)$-identity paths by means of the Curry-Howard mapping.

\begin{restatable}{corollary}{CCHId}
\label{CCHId} Let $s$ be a sort in $S$ and $\mathfrak{P}$ a path in $\mathrm{Pth}_{\boldsymbol{\mathcal{A}},s}$. The following statements are equivalent.
\begin{enumerate}
\item[(i)] $\mathfrak{P}$ is a $(1,0)$-identity path.
\item[(ii)] $\mathrm{CH}^{(1)}_{s}(\mathfrak{P})=\eta^{(1,0)\sharp}_{s}(\mathrm{sc}^{(0,1)}_{s}(\mathfrak{P}))$.
\item[(iii)] $\mathrm{CH}^{(1)}_{s}(\mathfrak{P})=\eta^{(1,0)\sharp}_{s}(\mathrm{tg}^{(0,1)}_{s}(\mathfrak{P}))$.
\item[(iv)] $\mathrm{CH}^{(1)}_{s}(\mathfrak{P})\in \eta^{(1,0)\sharp}[\mathrm{T}_{\Sigma}(X)]_{s}$.
\end{enumerate}
\end{restatable}

\begin{proof}
We first prove that (i) implies (ii). Let $\mathfrak{P}$ be a $(1,0)$-identity path of sort $s$ then, for some term $P\in\mathrm{T}_{\Sigma}(X)_{s}$, we have that $\mathfrak{P}=\mathrm{ip}^{(1,0)\sharp}_{s}(P)$.

The following chain of equalities holds
\begin{align*}
\mathrm{CH}^{(1)}_{s}
\left(
\mathfrak{P}
\right)&=
\mathrm{CH}^{(1)}_{s}\left(
\mathrm{ip}^{(1,0)\sharp}_{s}\left(
P
\right)\right)
\tag{1}
\\&=
\eta^{(1,0)\sharp}_{s}(P)
\tag{2}
\\&=
\eta^{(1,0)\sharp}_{s}
\left(
\mathrm{sc}^{(0,1)}_{s}\left(
\mathfrak{P}
\right)
\right).
\tag{3}
\end{align*}
The first equality unravels the description of $\mathfrak{P}$ as a $(1,0)$-identity path; the second equality follows from Proposition~\ref{PCHId}; finally, the last equality recovers the $(0,1)$-source of the path $\mathfrak{P}$.

That (ii) implies (iv) is straightforward, since $\mathrm{sc}^{(0,1)}_{s}(\mathfrak{P})$ is a term in $\mathrm{T}_{\Sigma}(X)_{s}$.

Finally, we prove that (iv) implies (i). By contraposition, let us assume that $\mathfrak{P}$ is a non-$(1,0)$-identity path. Then, by Proposition~\ref{PCHRewId}, 
$$\eta^{(1,\mathcal{A})}
\left[
\mathcal{A}
\right]\cap\mathrm{Subt}_{\Sigma^{\boldsymbol{\mathcal{A}}}}
\left(
\mathrm{CH}^{(1)}_{s}\left(
\mathfrak{P}
\right)\right)\neq\varnothing^{S},$$
contradicting the assumption that the term $\mathrm{CH}^{(1)}_{s}(\mathfrak{P})$ belongs to $\eta^{(1,0)\sharp}[\mathrm{T}_{\Sigma}(X)]_{s}$.

The proof that (i) implies (iii) and (iii) implies (iv) is as above. The remaining cases follow from the fact that the $(0,1)$-source and the $(0,1)$-target coincide for the $(1,0)$-identity paths.

This finishes the proof.
\end{proof}

We can now complete the characterization of the non-$(1,0)$-identity paths.

\begin{corollary}\label{CCHRew}
Let $s$ be a sort in $S$ and $\mathfrak{P}$ a path in $\mathrm{Pth}_{\boldsymbol{\mathcal{A}},s}$. The following statements are equivalent.
\begin{enumerate}
\item[(i)] $\mathfrak{P}$ is a  non-$(1,0)$-identity path.
\item[(ii)] $\eta^{(1,\mathcal{A})}[\mathcal{A}]\cap\mathrm{Subt}_{\Sigma^{\boldsymbol{\mathcal{A}}}}(\mathrm{CH}^{(1)}_{s}(\mathfrak{P}))\neq\varnothing^{S}$.
\end{enumerate}
\end{corollary}

\begin{proof}
That (i) implies (ii) was proved in Proposition~\ref{PCHRewId}.

Now we prove that (ii) implies (i). By contraposition, let us assume that $\mathfrak{P}$ is a $(1,0)$-identity path. Then, by Corollary~\ref{CCHId}, we have that $\mathrm{CH}^{(1)}_{s}(\mathfrak{P})$ is a term in $\eta^{(1,0)\sharp}[\mathrm{T}_{\Sigma}(X)]_{s}$, thus it cannot contain subterms of any sort from $\eta^{(1,\mathcal{A})}[\mathcal{A}]$.
\end{proof}

Another interesting result, that follows from Corollary~\ref{CCHId}, is the following: If the image under the Curry-Howard mapping of a path is the same as that of a $(1,0)$-identity path, then both paths are equal.

\begin{restatable}{corollary}{CCHUZId}
\label{CCHUZId}
Let $s$ be a sort in $S$ and $(\mathfrak{P},\mathfrak{Q})\in\mathrm{Ker}(\mathrm{CH}^{(1)})_{s}$. If $\mathfrak{P}$ or $\mathfrak{Q}$ is a $(1,0)$-identity path, then $\mathfrak{Q}=\mathfrak{P}$.
\end{restatable}
\begin{proof}
Let us assume that $\mathfrak{P}$ is a $(1,0)$-identity path. Then there exists a term $P$ in $\mathrm{T}_{\Sigma}(X)_{s}$ for which $\mathfrak{P}=\mathrm{ip}^{(1,0)\sharp}_{s}(P)$.

According to Corollary~\ref{CCHId}, we have that 
$$
\mathrm{CH}^{(1)}_{s}\left(
\mathfrak{P}
\right)
=
\eta^{(1,0)\sharp}_{s}\left(
P
\right).
$$

But, since $(\mathfrak{P},\mathfrak{Q})\in\mathrm{Ker}(\mathrm{CH}^{(1)})_{s}$, we have that $\mathrm{CH}^{(1)}_{s}(\mathfrak{Q})$ is a term in $\eta^{(1,0)\sharp}[\mathrm{T}_{\Sigma}(X)]_{s}$. Hence, by Corollary~\ref{CCHRew}, $\mathfrak{Q}$ is a $(1,0)$-identity path. Therefore, there exists a term $Q$ in $\mathrm{T}_{\Sigma}(X)_{s}$ for which $\mathfrak{Q}=\mathrm{ip}^{(1,0)\sharp}_{s}(Q)$.

In the same way, we have that 
$$
\mathrm{CH}^{(1)}_{s}\left(
\mathfrak{Q}
\right)
=
\eta^{(1,0)\sharp}_{s}\left(
Q
\right).
$$

Now, taking into account that $(\mathfrak{P},\mathfrak{Q})\in\mathrm{Ker}(\mathrm{CH}^{(1)})_{s}$ and the fact that, by Proposition~\ref{PEmb}, $\eta^{(1,0)\sharp}$ is injective, we have that 
$
P=Q.
$ 
Thus, we can affirm that 
$$
\mathfrak{P}
=
\mathrm{ip}^{(1,0)\sharp}_{s}(P)
=
\mathrm{ip}^{(1,0)\sharp}_{s}(Q)
=
\mathfrak{Q}.
$$

The statement follows.
\end{proof}

We will prove below that the non-$(1,0)$-identity paths containing echelons also have a nice characterization by means the Curry-Howard mapping. But, to do so, we first need to define some special  
$S$-sorted subsets of $\mathrm{T}_{\Sigma^{\boldsymbol{\mathcal{A}}}}(X)$ that are related to the binary operation of $0$-composition and to the $S$-sorted set of rewrite rules.

\begin{definition}\label{DA}
For the free $\Sigma^{\boldsymbol{\mathcal{A}}}$-algebra $\mathbf{T}_{\Sigma^{\boldsymbol{\mathcal{A}}}}(X)$
we define the following $S$-sorted subsets of $\mathrm{T}_{\Sigma^{\boldsymbol{\mathcal{A}}}}(X)$.

\textsf{(i)} The \emph{initial image of $\mathcal{A}$}, denoted by $\eta^{(1,\mathcal{A})}[\mathcal{A}]^{\mathrm{int}}$, is the $S$-sorted subset of $\mathrm{T}_{\Sigma^{\boldsymbol{\mathcal{A}}}}(X)$ defined, for every sort $s\in S$, as
$$
\eta^{(1,\mathcal{A})}[\mathcal{A}]^{\mathrm{int}}_{s}=
\left\lbrace
Q
\circ^{0\mathbf{T}_{\Sigma^{\boldsymbol{\mathcal{A}}}}(X)}_{s}
\mathfrak{p}^{\mathbf{T}_{\Sigma^{\boldsymbol{\mathcal{A}}}}(X)}
\,
\middle|
\,
\begin{gathered}
Q\in\mathrm{T}_{\Sigma^{\boldsymbol{\mathcal{A}}}}(X)_{s},
\\
\mathfrak{p}\in\mathcal{A}_{s}
\end{gathered}
\right\rbrace.
$$

\textsf{(ii)} The \emph{prescriptive image of $\mathcal{A}$}, denoted by $\eta^{(1,\mathcal{A})}[\mathcal{A}]^{\mathrm{pct}}$, is the $S$-sorted subset of $\mathrm{T}_{\Sigma^{\boldsymbol{\mathcal{A}}}}(X)$ defined, for every sort $s\in S$, as
$$
\eta^{(1,\mathcal{A})}[\mathcal{A}]^{\mathrm{pct}}_{s}=
\left\lbrace
\begin{gathered}
\mathfrak{p}^{\mathbf{T}_{\Sigma^{\boldsymbol{\mathcal{A}}}}(X)},
\\
Q
\circ^{0\mathbf{T}_{\Sigma^{\boldsymbol{\mathcal{A}}}}(X)}_{s}
\mathfrak{p}^{\mathbf{T}_{\Sigma^{\boldsymbol{\mathcal{A}}}}(X)},
\\
\left(R
\circ^{0\mathbf{T}_{\Sigma^{\boldsymbol{\mathcal{A}}}}(X)}_{s}
\mathfrak{p}^{\mathbf{T}_{\Sigma^{\boldsymbol{\mathcal{A}}}}(X)}
\right)
\circ^{0\mathbf{T}_{\Sigma^{\boldsymbol{\mathcal{A}}}}(X)}_{s}
Q,
\\
\mathfrak{p}^{\mathbf{T}_{\Sigma^{\boldsymbol{\mathcal{A}}}}(X)}
\circ^{0\mathbf{T}_{\Sigma^{\boldsymbol{\mathcal{A}}}}(X)}_{s}
Q
\end{gathered}
\,
\middle|
\,
\begin{gathered}
R,Q\in\mathrm{T}_{\Sigma^{\boldsymbol{\mathcal{A}}}}(X)_{s}, 
\\
\mathfrak{p}\in\mathcal{A}_{s}
\end{gathered}
\right\rbrace.
$$

\textsf{(iii)} The \emph{non-initial image of $\mathcal{A}$}, denoted by $\eta^{(1,\mathcal{A})}[\mathcal{A}]^{\neg\mathrm{int}}$, is the $S$-sorted subset of $\mathrm{T}_{\Sigma^{\boldsymbol{\mathcal{A}}}}(X)$ containing, for every sort $s\in S$, the  terms $P\in\mathrm{T}_{\Sigma^{\boldsymbol{\mathcal{A}}}}(X)_{s}$ of  the form
\begin{multline*}
\eta^{(1,\mathcal{A})}[\mathcal{A}]^{\neg\mathrm{int}}_{s}\\=
\left\lbrace
\begin{gathered}
\left(R
\circ^{0\mathbf{T}_{\Sigma^{\boldsymbol{\mathcal{A}}}}(X)}_{s}
\mathfrak{p}^{\mathbf{T}_{\Sigma^{\boldsymbol{\mathcal{A}}}}(X)}
\right)
\circ^{0\mathbf{T}_{\Sigma^{\boldsymbol{\mathcal{A}}}}(X)}_{s}
Q,
\\
\mathfrak{p}^{\mathbf{T}_{\Sigma^{\boldsymbol{\mathcal{A}}}}(X)}
\circ^{0\mathbf{T}_{\Sigma^{\boldsymbol{\mathcal{A}}}}(X)}_{s}
Q
\end{gathered}
\,
\middle|
\,
\begin{gathered}
R\in\mathrm{T}_{\Sigma^{\boldsymbol{\mathcal{A}}}}(X)_{s},
\\
Q\in \mathrm{T}_{\Sigma^{\boldsymbol{\mathcal{A}}}}(X)_{s}-\eta^{(1,\mathcal{A})}[\mathcal{A}]^{\mathrm{pct}}_{s},
\\
\mathfrak{p}\in\mathcal{A}_{s}
\end{gathered}
\right\rbrace.
\end{multline*}
\end{definition}

To prove the following lemmas it, essentially, suffices to unravel the definition of the Curry-Howard mapping at the suitable path. Thus, we leave the proof to the reader.

\begin{lemma}\label{LCHEch}
Let $s$ be a sort in $S$ and $\mathfrak{P}$ a path in $\mathrm{Pth}_{\boldsymbol{\mathcal{A}},s}$. The following statements are equivalent.
\begin{enumerate}
\item[(i)] $\mathfrak{P}$ contains echelons.
\item[(ii)] $\mathrm{CH}^{(1)}_{s}(\mathfrak{P})\in \eta^{(1,\mathcal{A})}[\mathcal{A}]^{\mathrm{pct}}_{s}$.
\end{enumerate}
\end{lemma}

\begin{lemma}\label{LCHUEch}
Let $s$ be a sort in $S$ and $\mathfrak{P}$ a path in $\mathrm{Pth}_{\boldsymbol{\mathcal{A}},s}$. The following statements are equivalent.
\begin{enumerate}
\item[(i)] $\mathfrak{P}$ is an echelon.
\item[(ii)] $\mathrm{CH}^{(1)}_{s}(\mathfrak{P})\in \eta^{(1,\mathcal{A})}[\mathcal{A}]_{s}$.
\end{enumerate}
\end{lemma}

\begin{lemma}\label{LCHEchInt}
Let $s$ be a sort in $S$ and $\mathfrak{P}$ a path in $\mathrm{Pth}_{\boldsymbol{\mathcal{A}},s}$. The following statements are equivalent.
\begin{enumerate}
\item[(i)] $\bb{\mathfrak{P}}>1$ and $\mathfrak{P}$ has its first echelon on its first step.
\item[(ii)] $\mathrm{CH}^{(1)}_{s}(\mathfrak{P})\in \eta^{(1,\mathcal{A})}[\mathcal{A}]^{\mathrm{int}}_{s}$.
\end{enumerate}
\end{lemma}

\begin{lemma}\label{LCHEchNInt}
Let $s$ be a sort in $S$ and $\mathfrak{P}$ a path in $\mathrm{Pth}_{\boldsymbol{\mathcal{A}},s}$. The following statements are equivalent.
\begin{enumerate}
\item[(i)] $\bb{\mathfrak{P}}>1$ and $\mathfrak{P}$ has its first echelon in a different step from the initial one.
\item[(ii)] $\mathrm{CH}^{(1)}_{s}(\mathfrak{P})\in \eta^{(1,\mathcal{A})}[\mathcal{A}]^{\neg\mathrm{int}}_{s}$.
\end{enumerate}
\end{lemma}

We will prove below that the echelonless paths also have a nice characterization by means of the Curry-Howard mapping. 

\begin{definition}\label{DAII}
Let $s$ be a sort in $S$, $\mathbf{s}$ a word in $S^{\star}-\{\lambda\}$, and $\sigma\in\Sigma^{\boldsymbol{\mathcal{A}}}_{\mathbf{s}, s}$. Then we let $\mathcal{T}\left(\sigma,\mathrm{T}_{\Sigma^{\boldsymbol{\mathcal{A}}}}(X)\right)_{1}$ stand for the subset
$$
\left\lbrace
\sigma^{\mathbf{T}_{\Sigma^{\boldsymbol{\mathcal{A}}}}(X)}
\left(\left(P_{j}
\right)_{j\in\bb{\mathbf{s}}}
\right)
\, 
\middle|
\,
\begin{gathered}
(P_{j})_{j\in\bb{\mathbf{s}}}
\in\mathrm{T}_{\Sigma^{\boldsymbol{\mathcal{A}}}}(X)_{\mathbf{s}},
\\
\exists j\in\bb{\mathbf{s}}\,
(P_{j}\not\in \eta^{(1,0)\sharp}
[\mathrm{T}_{\Sigma}(X)]_{s_{j}})
\end{gathered}
\right\rbrace
$$
of $\mathrm{T}_{\Sigma^{\boldsymbol{\mathcal{A}}}}(X)_{s}$.

From $\mathcal{T}\left(\sigma,\mathrm{T}_{\Sigma^{\boldsymbol{\mathcal{A}}}}(X)\right)_{1}$ we define the following $S$-sorted subset of $\mathrm{T}_{\Sigma^{\boldsymbol{\mathcal{A}}}}(X)$.

The \emph{echelonless $S$-sorted subset of $\mathrm{T}_{\Sigma^{\boldsymbol{\mathcal{A}}}}(X)$}, written  $\mathrm{T}_{\Sigma^{\boldsymbol{\mathcal{A}}}}(X)^{
\mathsf{E}
}$, is  defined, for every sort $s\in S$, as
$$
\mathrm{T}_{\Sigma^{\boldsymbol{\mathcal{A}}}}(X)
^{\mathsf{E}}_{s}=
\bigcup_{\sigma\in\Sigma_{\neq\lambda, s}}
\mathcal{T}
\left(
\sigma,\mathrm{T}_{\Sigma^{\boldsymbol{\mathcal{A}}}}(X)
\right)_{1}.
$$
\end{definition}

To prove the following lemma it, essentially, suffices to unravel the definition of the Curry-Howard mapping at the suitable path. Thus, we leave the proof to the reader.

\begin{lemma}\label{LCHNEch}
Let $s$ be a sort in $S$ and $\mathfrak{P}$ a path in $\mathrm{Pth}_{\boldsymbol{\mathcal{A}},s}$. The following statements are equivalent.
\begin{enumerate}
\item[(i)] $\mathfrak{P}$ is an echelonless  path.
\item[(ii)] $\mathrm{CH}^{(1)}_{s}(\mathfrak{P})\in \mathrm{T}_{\Sigma^{\boldsymbol{\mathcal{A}}}}(X)^{\mathsf{E}}_{s}
$.
\end{enumerate}
\end{lemma}

Another interesting result that follows from the just stated lemmas is that the Curry-Howard mapping
restricted to $\mathrm{Ech}^{(1,\mathcal{A})}[\mathcal{A}]$, the subset of $\mathrm{Pth}_{\boldsymbol{\mathcal{A}}}$ formed by the echelons, is an injective mapping.

\begin{restatable}{proposition}{PCHEch}
\label{PCHEch}
Let $s$ be a sort in $S$ and $(\mathfrak{P},\mathfrak{Q})\in\mathrm{Ker}(\mathrm{CH}^{(1)})_{s}$. If 
$\mathfrak{P}$ or $\mathfrak{Q}$ is an echelon, then $\mathfrak{Q}=\mathfrak{P}$.
\end{restatable}
\begin{proof}
Assume that $\mathfrak{P}$ is an echelon, i.e., that $\mathfrak{P}=\mathrm{ech}^{(1,\mathcal{A})}_{s}(\mathfrak{p})$, for some  rewrite rule $\mathfrak{p}\in\mathcal{A}_{s}$. Then, by Definition~\ref{DCH}, we have that 
$$
\mathrm{CH}^{(1)}_{s}\left(\mathfrak{P}\right)=
\mathfrak{p}^{\mathbf{T}_{\Sigma^{\boldsymbol{\mathcal{A}}}}(X)}.
$$

Since $(\mathfrak{P},\mathfrak{Q})\in\mathrm{Ker}(\mathrm{CH}^{(1)})_{s}$, we have that $\mathrm{CH}^{(1)}_{s}(\mathfrak{Q})$ is a term in $\eta^{(1,\mathcal{A})}[\mathcal{A}]_{s}$. Hence, by Lemma~\ref{LCHUEch}, $\mathfrak{Q}$ is an echelon. Since $\mathrm{CH}^{(1)}_{s}(\mathfrak{P})=\mathrm{CH}^{(1)}_{s}(\mathfrak{Q})$, we can affirm that the unique rewrite rule appearing in both $\mathfrak{P}$ and $\mathfrak{Q}$ is the same. Consequently, $\mathfrak{P}=\mathfrak{Q}$.
\end{proof}

Before proving, in the following section, that $\mathrm{Ker}(\mathrm{CH}^{(1)})$ is a closed $\Sigma^{\boldsymbol{\mathcal{A}}}$-congruence on the partial $\Sigma^{\boldsymbol{\mathcal{A}}}$-algebra 
$\mathbf{Pth}_{\boldsymbol{\mathcal{A}}}$, we need to prove, in the following lemma, that whenever two paths of the same sort are in the kernel of the Curry-Howard mapping, then they have the same length, the same $(0,1)$-source, and the same $(0,1)$-target.

\begin{restatable}{lemma}{LCH}
\label{LCH}
Let $s$ a sort in $S$ and $(\mathfrak{P},\mathfrak{Q})\in\mathrm{Ker}(\mathrm{CH}^{(1)})_{s}$. Then the following statements hold.
\begin{enumerate}
\item[(i)] $\bb{\mathfrak{P}}=\bb{\mathfrak{Q}}$;
\item[(ii)] $\mathrm{sc}^{(0,1)}_{s}(\mathfrak{P})=\mathrm{sc}^{(0,1)}_{s}(\mathfrak{Q})$;
\item[(iii)] $\mathrm{tg}^{(0,1)}_{s}(\mathfrak{P})=\mathrm{tg}^{(0,1)}_{s}(\mathfrak{Q})$.
\end{enumerate}
\end{restatable}

\begin{proof}
If either $\mathfrak{P}$ or $\mathfrak{Q}$ is a $(1,0)$-identity path, then the statement follows according to Corollary~\ref{CCHUZId}. Therefore, we can assume that neither $\mathfrak{P}$ nor $\mathfrak{Q}$ are $(1,0)$-identity paths.

We prove the statements by Artinian induction on $(\coprod\mathrm{Pth}_{\boldsymbol{\mathcal{A}}},\leq_{\mathbf{Pth}_{\boldsymbol{\mathcal{A}}}})$.

\textsf{Base step of the Artinian induction}.

Let $(\mathfrak{P},s)$ be a minimal element of $(\coprod\mathrm{Pth}_{\boldsymbol{\mathcal{A}}}, \leq_{\mathbf{Pth}_{\boldsymbol{\mathcal{A}}}})$. Then, by Proposition~\ref{PMinimal} and taking into account that $\mathfrak{P}$ is a non-$(1,0)$-identity path, we have that $\mathfrak{P}$ is an echelon. The statement follows in virtue of Proposition~\ref{PCHEch}.

This completes the base step.

\textsf{Inductive step of the Artinian induction}.

Let $(\mathfrak{P},s)$ be a non-minimal element in $(\coprod\mathrm{Pth}_{\boldsymbol{\mathcal{A}}}, \leq_{\mathbf{Pth}_{\boldsymbol{\mathcal{A}}}})$. Let us suppose that, for every sort $t\in S$ and every path $\mathfrak{P}'\in\mathrm{Pth}_{\boldsymbol{\mathcal{A}},t}$, if $(\mathfrak{P}',t)<_{\mathbf{Pth}_{\boldsymbol{\mathcal{A}}}}(\mathfrak{P},s)$, then the statement holds for $\mathfrak{P}'$, i.e., for every path $\mathfrak{Q}'$ in $\mathrm{Pth}_{\boldsymbol{\mathcal{A}},t}$ if $(\mathfrak{P}',\mathfrak{Q}')\in\mathrm{Ker}(\mathrm{CH}^{(1)})_{t}$, then $\mathfrak{P}'$ and $\mathfrak{Q}'$ have the same length, the same $(0,1)$-source and the same $(0,1)$-target.

Since $(\mathfrak{P},s)$ is a non-minimal element of $(\coprod\mathrm{Pth}_{\boldsymbol{\mathcal{A}}}, \leq_{\mathbf{Pth}_{\boldsymbol{\mathcal{A}}}})$ and taking into account that $\mathfrak{P}$ is a non-$(1,0)$-identity path, we have, by Lemma~\ref{LOrdI}, that $\mathfrak{P}$ is either~(1) a path of length strictly greater than one containing at least one echelon or~(2) an echelonless path.

If~(1), then let $i\in \bb{\mathfrak{P}}$ be the first index for which the one-step subpath $\mathfrak{P}^{i,i}$ is an echelon. We distinguish two cases accordingly.

If $i=0$, i.e., if $\mathfrak{P}$ has its first echelon on its first step, then, according to Definition~\ref{DCH}, we have that
$$
\mathrm{CH}^{(1)}_{s}\left(
\mathfrak{P}
\right)=
\mathrm{CH}^{(1)}_{s}\left(
\mathfrak{P}^{1,\bb{\mathfrak{P}}-1}
\right)
\circ_{s}^{0\mathbf{T}_{\Sigma^{\boldsymbol{\mathcal{A}}}}(X)}
\mathrm{CH}^{(1)}_{s}\left(
\mathfrak{P}^{0,0}\right).
$$

Since $\mathrm{CH}^{(1)}_{s}(\mathfrak{P})\in\eta^{(1,\mathcal{A})}[\mathcal{A}]^{\mathrm{int}}_{s}$ and $(\mathfrak{P},\mathfrak{Q})\in\mathrm{Ker}(\mathrm{CH}^{(1)})_{s}$, we have, by Lemma~\ref{LCHEchInt} that $\mathfrak{Q}$ is a path of length strictly greater than one containing an echelon on its first step. 

Thus, according to Definition~\ref{DCH}, we have that 
$$
\mathrm{CH}^{(1)}_{s}\left(
\mathfrak{Q}
\right)=
\mathrm{CH}^{(1)}_{s}\left(
\mathfrak{Q}^{1,\bb{\mathfrak{Q}}-1}
\right)
\circ_{s}^{0\mathbf{T}_{\Sigma^{\boldsymbol{\mathcal{A}}}}(X)}
\mathrm{CH}^{(1)}_{s}\left(
\mathfrak{Q}^{0,0}
\right).
$$

Since $(\mathfrak{P},\mathfrak{Q})\in\mathrm{Ker}(\mathrm{CH}^{(1)})_{s}$, we have that $(\mathfrak{P}^{1,\bb{\mathfrak{P}}-1},\mathfrak{Q}^{1,\bb{\mathfrak{Q}}-1})$ and $(\mathfrak{P}^{0,0},\mathfrak{Q}^{0,0})$ are in $\mathrm{Ker}(\mathrm{CH}^{(1)})_{s}$. Note that, according to Definition~\ref{DOrd}, we have that $(\mathfrak{P}^{0,0},s)$ and $(\mathfrak{P}^{1,\bb{\mathfrak{P}}-1},s)$ $\prec_{\mathbf{Pth}_{\boldsymbol{\mathcal{A}}}}$-precede $(\mathfrak{P},s)$.

Therefore, by the inductive hypothesis, the paths $\mathfrak{P}^{0,0}$ and $\mathfrak{Q}^{0,0}$, and the paths $\mathfrak{P}^{1,\bb{\mathfrak{P}}-1}$ and $\mathfrak{Q}^{1,\bb{\mathfrak{Q}}-1}$ have, respectively, the same length, the same $(0,1)$-source and the same $(0,1)$-target.

Therefore, we have that
\begin{itemize}
\item[(i)] $\bb{\mathfrak{P}}=
\bb{\mathfrak{P}^{0,0}}+
\bb{\mathfrak{P}^{1,\bb{\mathfrak{P}}-1}}
=
\bb{\mathfrak{Q}^{0,0}}+
\bb{\mathfrak{Q}^{1,\bb{\mathfrak{Q}}-1}}
=\bb{\mathfrak{Q}};$
\item[(ii)] $\mathrm{sc}^{(0,1)}_{s}\left(\mathfrak{P}
\right)
=
\mathrm{sc}^{(0,1)}_{s}\left(
\mathfrak{P}^{0,0}
\right)
=
\mathrm{sc}^{(0,1)}_{s}\left(
\mathfrak{Q}^{0,0}
\right)
=
\mathrm{sc}^{(0,1)}_{s}\left(
\mathfrak{Q}
\right);$
\item[(iii)] $\mathrm{tg}^{(0,1)}_{s}\left(
\mathfrak{P}\right)
=
\mathrm{tg}^{(0,1)}_{s}\left(
\mathfrak{P}^{1,\bb{\mathfrak{P}}-1}
\right)
=
\mathrm{tg}^{(0,1)}_{s}\left(
\mathfrak{Q}^{1,\bb{\mathfrak{Q}}-1}
\right)
=
\mathrm{tg}^{(0,1)}_{s}\left(
\mathfrak{Q}\right).$
\end{itemize}

The case of $\mathfrak{P}$ being a  path of length strictly greater than one containing its first echelon on its first step follows.

If $i>0$, that is, if $\mathfrak{P}$ is a  path of length strictly greater than one containing its first  echelon on a step different from the initial one, then, according to Definition~\ref{DCH}, we have that
$$
\mathrm{CH}^{(1)}_{s}\left(
\mathfrak{P}
\right)=
\mathrm{CH}^{(1)}_{s}\left(
\mathfrak{P}^{i,\bb{\mathfrak{P}}-1}
\right)
\circ_{s}^{0\mathbf{T}_{\Sigma^{\boldsymbol{\mathcal{A}}}}(X)}
\mathrm{CH}^{(1)}_{s}\left(
\mathfrak{P}^{0,i-1}
\right).
$$

Since $\mathrm{CH}^{(1)}_{s}(\mathfrak{P})\in\eta^{(1,\mathcal{A})}[\mathcal{A}]^{\neg\mathrm{int}}_{s}$ and $(\mathfrak{P},\mathfrak{Q})\in\mathrm{Ker}(\mathrm{CH}^{(1)})_{s}$, we have, by Lemma~\ref{LCHEchNInt}, that $\mathfrak{Q}$ is a path of length strictly greater than one containing its first echelon on a step different from the initial one.

Thus, if $j\in\bb{\mathfrak{Q}}$ is the first index for which $\mathfrak{Q}^{j,j}$ is an  echelon, according to Definition~\ref{DCH}, we have that
$$
\mathrm{CH}^{(1)}_{s}\left(
\mathfrak{Q}
\right)
=
\mathrm{CH}^{(1)}_{s}\left(
\mathfrak{Q}^{j,\bb{\mathfrak{Q}}-1}
\right)
\circ_{s}^{0\mathbf{T}_{\Sigma^{\boldsymbol{\mathcal{A}}}}(X)}
\mathrm{CH}^{(1)}_{s}\left(
\mathfrak{Q}^{0,j-1}
\right).
$$

Since $(\mathfrak{P},\mathfrak{Q})\in\mathrm{Ker}(\mathrm{CH}^{(1)})_{s}$, we have that $(\mathfrak{P}^{i,\bb{\mathfrak{P}}-1},\mathfrak{Q}^{j,\bb{\mathfrak{Q}}-1})$ and $(\mathfrak{P}^{0,i-1}, \mathfrak{Q}^{0,j-1})$ are in $\mathrm{Ker}(\mathrm{CH}^{(1)})_{s}$. Note that, according to Definition~\ref{DOrd}, we have that $(\mathfrak{P}^{0,i-1},s)$ and $(\mathfrak{P}^{i,\bb{\mathfrak{P}}-1},s)$ $\prec_{\mathbf{Pth}_{\boldsymbol{\mathcal{A}}}}$-precede $(\mathfrak{P},s)$. 

Therefore, by the inductive hypothesis, the  paths $\mathfrak{P}^{0,i-1}$ and $\mathfrak{Q}^{0,j-1}$, and the  paths $\mathfrak{P}^{i,\bb{\mathfrak{P}}-1}$ and $\mathfrak{Q}^{j,\bb{\mathfrak{Q}}-1}$ have, respectively, the same length, the same $(0,1)$-source and the same $(0,1)$-target. In particular $i=j$.

Therefore, we have that
\begin{itemize}
\item[(i)] $\bb{\mathfrak{P}}=
\bb{\mathfrak{P}^{0,i-1}}+
\bb{\mathfrak{P}^{i,\bb{\mathfrak{P}}-1}}
=
\bb{\mathfrak{Q}^{0,j-1}}+
\bb{\mathfrak{Q}^{j,\bb{\mathfrak{Q}}-1}}
=\bb{\mathfrak{Q}};$
\item[(ii)] $\mathrm{sc}^{(0,1)}_{s}\left(\mathfrak{P}
\right)
=
\mathrm{sc}^{(0,1)}_{s}\left(
\mathfrak{P}^{0,i-1}
\right)
=
\mathrm{sc}^{(0,1)}_{s}\left(
\mathfrak{Q}^{0,j-1}
\right)
=
\mathrm{sc}^{(0,1)}_{s}\left(
\mathfrak{Q}
\right);$
\item[(iii)] $\mathrm{tg}^{(0,1)}_{s}\left(
\mathfrak{P}\right)
=
\mathrm{tg}^{(0,1)}_{s}\left(
\mathfrak{P}^{i,\bb{\mathfrak{P}}-1}
\right)
=
\mathrm{tg}^{(0,1)}_{s}\left(
\mathfrak{Q}^{j,\bb{\mathfrak{Q}}-1}
\right)
=
\mathrm{tg}^{(0,1)}_{s}\left(
\mathfrak{Q}\right).$
\end{itemize}

The case of $\mathfrak{P}$ being a path of length strictly greater than one containing its first echelon on a step different from the initial one follows.

Case~(1) follows.

If~(2), then in virtue of Lemma~\ref{LPthHeadCt} there exists a unique word $\mathbf{s}\in S^{\star}-\{\lambda\}$ and a unique operation symbol $\sigma\in \Sigma_{\mathbf{s},s}$ associated to $\mathfrak{P}$. Let $(\mathfrak{P}_{j})_{j\in\bb{\mathbf{s}}}$ be the family of paths in $\mathrm{Pth}_{\boldsymbol{\mathcal{A}},\mathbf{s}}$ which, in virtue of Lemma~\ref{LPthExtract}, we can extract from $\mathfrak{P}$. Then, according to Definition~\ref{DCH}, the value of the Curry-Howard mapping at $\mathfrak{P}$ is given by
$$
\mathrm{CH}^{(1)}_{s}\left(
\mathfrak{P}
\right)=
\sigma^{\mathbf{T}_{\Sigma^{\boldsymbol{\mathcal{A}}}}(X)}
\left(\left(\mathrm{CH}^{(1)}_{s_{j}}\left(
\mathfrak{P}_{j}
\right)\right)_{j\in\bb{\mathbf{s}}}
\right).
$$

Since $\mathrm{CH}^{(1)}_{s}(\mathfrak{P})\in\mathcal{T}(\sigma,\mathrm{T}_{\Sigma^{\boldsymbol{\mathcal{A}}}}(X))_{1}$, which is a subset of $\mathrm{T}_{\Sigma^{\boldsymbol{\mathcal{A}}}}(X)^{\mathsf{E}}_{s}$, and $(\mathfrak{P}, \mathfrak{Q})\in\mathrm{Ker}(\mathrm{CH}^{(1)})_{s}$ we have, by Lemma~\ref{LCHNEch}, that $\mathfrak{Q}$ is an  echelonless path associated to $\sigma$, the same operation symbol as that associated to $\mathfrak{P}$.  

Let $(\mathfrak{Q}_{j})_{j\in\bb{\mathbf{s}}}$ be the family of  paths in $\mathrm{Pth}_{\boldsymbol{\mathcal{A}},\mathbf{s}}$ which, by Lemma~\ref{LPthExtract}, we can extract from $\mathfrak{Q}$. According to Definition~\ref{DCH}, the image of the Curry-Howard mapping at $\mathfrak{Q}$ is given by 
$$
\mathrm{CH}^{(1)}_{s}\left(
\mathfrak{Q}
\right)
=
\sigma^{\mathbf{T}_{\Sigma^{\boldsymbol{\mathcal{A}}}}(X)}
\left(\left(
\mathrm{CH}^{(1)}_{s_{j}}\left(
\mathfrak{Q}_{j}
\right)\right)_{j\in\bb{\mathbf{s}}}
\right).
$$

Since $(\mathfrak{P},\mathfrak{Q})\in\mathrm{Ker}(\mathrm{CH}^{(1)})_{s}$, we have, for every $j\in\bb{\mathbf{s}}$, that $(\mathfrak{P}_{j}, \mathfrak{Q}_{j})\in\mathrm{Ker}(\mathrm{CH}^{(1)})_{s_{j}}$. Note that, according to Definition~\ref{DOrd}, we have that, for every $j\in\bb{\mathbf{s}}$, $(\mathfrak{P}_{j}, s_{j})$ $\prec_{\mathbf{Pth}_{\boldsymbol{\mathcal{A}}}}$-precedes $(\mathfrak{P},s)$.

Therefore, by the inductive hypothesis, for every $j\in\bb{\mathbf{s}}$, the  paths $\mathfrak{P}_{j}$ and $\mathfrak{Q}_{j}$ have the same length, the same $(0,1)$-source and the same $(0,1)$-target.

Therefore we have that
\leqnomode
\allowdisplaybreaks
\begin{align*}
\bb{\mathfrak{P}}
&=
\sum_{j\in\bb{\mathbf{s}}}\bb{\mathfrak{P}_{j}}
=\sum_{j\in\bb{\mathbf{s}}}\bb{\mathfrak{Q}_{j}}
=\bb{\mathfrak{Q}};
\tag{i}
\\
\mathrm{sc}^{(0,1)}_{s}\left(
\mathfrak{P}
\right)&=
\sigma^{\mathbf{T}_{\Sigma}(X)}
\left(\left(\mathrm{sc}^{(0,1)}_{s_{j}}
\left(
\mathfrak{P}_{j}
\right)\right)_{j\in\bb{\mathbf{s}}}
\right)
\tag{ii}
\\&=\sigma^{\mathbf{T}_{\Sigma}(X)}
\left(\left(\mathrm{sc}^{(0,1)}_{s_{j}}\left(
\mathfrak{Q}_{j}
\right)\right)_{j\in\bb{\mathbf{s}}}
\right)
\\&=
\mathrm{sc}^{(0,1)}_{s}\left(
\mathfrak{Q}\right)
;
\\
\mathrm{tg}^{(0,1)}_{s}\left(
\mathfrak{P}
\right)&=
\sigma^{\mathbf{T}_{\Sigma}(X)}
\left(\left(\mathrm{tg}^{(0,1)}_{s_{j}}
\left(
\mathfrak{P}_{j}
\right)\right)_{j\in\bb{\mathbf{s}}}
\right)
\tag{iii}
\\&=\sigma^{\mathbf{T}_{\Sigma}(X)}
\left(\left(\mathrm{tg}^{(0,1)}_{s_{j}}\left(
\mathfrak{Q}_{j}
\right)\right)_{j\in\bb{\mathbf{s}}}
\right)
\\&=
\mathrm{tg}^{(0,1)}_{s}\left(
\mathfrak{Q}\right).
\end{align*}

The case of $\mathfrak{P}$ being an echelonless path  follows.

This finishes the proof.
\end{proof}

In the following proposition we prove that the Curry-Howard mapping restricted to the subset of 
$\mathrm{Pth}_{\boldsymbol{\mathcal{A}}}$ formed by the one-step paths is an injective mapping.

\begin{restatable}{proposition}{PCHOneStep}
\label{PCHOneStep}
Let $s$ be a sort in $S$ and $(\mathfrak{P},\mathfrak{Q})\in\mathrm{Ker}(\mathrm{CH}^{(1)})_{s}$. If $\mathfrak{P}$ or $\mathfrak{Q}$ is a one-step path, then $\mathfrak{Q}=\mathfrak{P}$.
\end{restatable}

\begin{proof}
Let $\mathfrak{P}$ be a one-step path.
We prove the statement by Artinian induction on $(\coprod\mathrm{Pth}_{\boldsymbol{\mathcal{A}}},\leq_{\mathbf{Pth}_{\boldsymbol{\mathcal{A}}}})$.

\textsf{Base step of the Artinian induction}.

Let $(\mathfrak{P},s)$ be a minimal element of $(\coprod\mathrm{Pth}_{\boldsymbol{\mathcal{A}}}, \leq_{\mathbf{Pth}_{\boldsymbol{\mathcal{A}}}})$. Then, by Proposition~\ref{PMinimal}, the path $\mathfrak{P}$ is either~(1) a $(1,0)$-identity path on a simple term, or~(2) an echelon. But, since we are assuming that $\mathfrak{P}$ is a one-step path, it follows that $\mathfrak{P}$ can only be an echelon. This case follows from Proposition~\ref{PCHEch}.

\textsf{Inductive step of the Artinian induction}.

Let $(\mathfrak{P},s)$ be a non-minimal element of $(\coprod\mathrm{Pth}_{\boldsymbol{\mathcal{A}}}, \leq_{\mathbf{Pth}_{\boldsymbol{\mathcal{A}}}})$. Let us suppose that, for every sort $t\in S$ and every path $\mathfrak{P}'\in\mathrm{Pth}_{\boldsymbol{\mathcal{A}},t}$, if $(\mathfrak{P}',t)<_{\mathbf{Pth}_{\boldsymbol{\mathcal{A}}}}(\mathfrak{P},s)$, then the statement holds for $\mathfrak{P}'$, i.e., if $\mathfrak{P}'$ is a one-step path in $\mathrm{Pth}_{\boldsymbol{\mathcal{A}},t}$, and $\mathfrak{Q}'$ is a path in $\mathrm{Pth}_{\boldsymbol{\mathcal{A}},t}$ such that $(\mathfrak{P}',\mathfrak{Q}')\in\mathrm{Ker}(\mathrm{CH}^{(1)})_{t}$, then $\mathfrak{P}'=\mathfrak{Q}'$.

Since $(\mathfrak{P},s)$ is a non-minimal element in $(\coprod\mathrm{Pth}_{\boldsymbol{\mathcal{A}}}, \leq_{\mathbf{Pth}_{\boldsymbol{\mathcal{A}}}})$, it follows, by Lemma~\ref{LOrdI}, that $\mathfrak{P}$ is either~(1) a $(1,0)$-identity path on a non-minimal term, or~(2) a path of length strictly greater than one containing at least one echelon or~(3) an echelonless path. But, since we are assuming that $\mathfrak{P}$ is a one-step path, it follows that $\mathfrak{P}$ can only be an echelonless path. 

Thus, following Lemma~\ref{LPthHeadCt}, there exists a unique word $\mathbf{s}\in S^{\star}-\{\lambda\}$ and a unique operation symbol $\sigma\in \Sigma_{\mathbf{s},s}$ associated to $\mathfrak{P}$. Let $(\mathfrak{P}_{j})_{j\in\bb{\mathbf{s}}}$ be the family of paths in $\mathrm{Pth}_{\boldsymbol{\mathcal{A}},\mathbf{s}}$ which, in virtue of Lemma~\ref{LPthExtract}, we can extract from $\mathfrak{P}$. In this case, the value of the Curry-Howard mapping at $\mathfrak{P}$ is given by 
$$
\mathrm{CH}^{(1)}_{s}\left(
\mathfrak{P}
\right)=
\sigma^{\mathbf{T}_{\Sigma^{\boldsymbol{\mathcal{A}}}}(X)}
\left(\left(\mathrm{CH}^{(1)}_{s_{j}}\left(
\mathfrak{P}_{j}
\right)\right)_{j\in\bb{\mathbf{s}}}
\right).
$$

Since $\mathfrak{P}$ is a one-step path, there exists a unique $k\in\bb{\mathbf{s}}$ for which  $\mathfrak{P}_{k}$ is also a one-step path. Moreover, for every $j\in\bb{\mathbf{s}}-\{k\}$, $\mathfrak{P}_{j}$ is a $(1,0)$-identity path. Let us note that, according to Definition~\ref{DOrd}, $(\mathfrak{P}_{k},s_{k})\prec_{\mathbf{Pth}_{\boldsymbol{\mathcal{A}}}}(\mathfrak{P},s)$.

Since $(\mathfrak{P},\mathfrak{Q})\in\mathrm{Ker}(\mathrm{CH}^{(1)})_{s}$, and $\mathrm{CH}^{(1)}_{s}(\mathfrak{P})\in \mathcal{T}(\sigma, \mathrm{T}_{\Sigma^{\boldsymbol{\mathcal{A}}}}(X))_{1}$, we conclude, in virtue of Lemma~\ref{LCHNEch}, that  $\mathfrak{Q}$ is an echelonless path associated to $\sigma$. In fact $\mathfrak{Q}$ is also a one-step path in virtue of Lemma~\ref{LCH}. Let $(\mathfrak{Q}_{j})_{j\in\bb{\mathbf{s}}}$ be the family of paths which we can extract from $\mathfrak{Q}$ in virtue of Lemma~\ref{LPthExtract}. By construction, the value of the Curry-Howard mapping at $\mathfrak{Q}$ is given by 
$$
\mathrm{CH}^{(1)}_{s}\left(
\mathfrak{Q}
\right)=
\sigma^{\mathbf{T}_{\Sigma^{\boldsymbol{\mathcal{A}}}}(X)}
\left(\left(\mathrm{CH}^{(1)}_{s_{j}}\left(
\mathfrak{Q}_{j}
\right)\right)_{j\in\bb{\mathbf{s}}}
\right).
$$

Thus, for every $j\in\bb{\mathbf{s}}$, the pairs $(\mathfrak{P}_{j},\mathfrak{Q}_{j})$ are in $\mathrm{Ker}(\mathrm{CH}^{(1)})_{s_{j}}$. Then, by Corollary~\ref{CCHUZId}, we have that, for every $j\in\bb{\mathbf{s}}-\{k\}$, $\mathfrak{P}_{j}=\mathfrak{Q}_{j}$. In particular, for the index $k\in\bb{\mathbf{s}}$ we have, by the inductive hypothesis, that $\mathfrak{P}_{k}=\mathfrak{Q}_{k}$.

Since $\mathfrak{P}$ and $\mathfrak{Q}$ are one-step paths, we have, by Corollary~\ref{CUStep}, that
$$
\mathfrak{P}=
\sigma^{\mathbf{Pth}_{\boldsymbol{\mathcal{A}}}}
\left(\left(\mathfrak{P}_{j}
\right)_{j\in\bb{\mathbf{s}}}
\right)=
\sigma^{\mathbf{Pth}_{\boldsymbol{\mathcal{A}}}}
\left(\left(\mathfrak{Q}_{j}
\right)
_{j\in\bb{\mathbf{s}}}
\right)
=\mathfrak{Q}.
$$

This finishes the proof.
\end{proof}

\section{A congruence on paths}
We are now in position to prove that the kernel of the Curry-Howard mapping is a closed 
$\Sigma^{\boldsymbol{\mathcal{A}}}$-congruence on $\mathbf{Pth}_{\boldsymbol{\mathcal{A}}}$.

\begin{restatable}{proposition}{PCHCong}
\label{PCHCong}
$\mathrm{Ker}(\mathrm{CH}^{(1)})$ is a closed $\Sigma^{\boldsymbol{\mathcal{A}}}$-congruence on the partial 
$\Sigma^{\boldsymbol{\mathcal{A}}}$-algebra $\mathbf{Pth}_{\boldsymbol{\mathcal{A}}}$.
\end{restatable}

\begin{proof}
By Proposition~\ref{PCHHom}, the Curry-Howard mapping is a $\Sigma$-homomorphism from $\mathbf{Pth}^{(0,1)}_{\boldsymbol{\mathcal{A}}}$, the $\Sigma$-reduct of the $\Sigma^{\boldsymbol{\mathcal{A}}}$-algebra $\mathbf{Pth}_{\boldsymbol{\mathcal{A}}}$, to $\mathbf{T}_{\Sigma^{\boldsymbol{\mathcal{A}}}}^{(0,1)}(X)$, the $\Sigma$-reduct of the free $\Sigma^{\boldsymbol{\mathcal{A}}}$-algebra $\mathbf{T}_{\Sigma^{\boldsymbol{\mathcal{A}}}}(X)$. Therefore, $\mathrm{Ker}(\mathrm{CH}^{(1)})$ is a $\Sigma$-congruence.

Thus, we only need to check the compatibility of $\mathrm{Ker}(\mathrm{CH}^{(1)})$ with the categorial operations. The addition of new constants to the signature, obviously, does not matter. 

Regarding the $0$-source and the $0$-target operations, we note that, for every sort $s\in S$, if $(\mathfrak{P},\mathfrak{Q})$ is a pair of paths in $\mathrm{Ker}(\mathrm{CH}^{(1)})_{s}$ then, by Lemma~\ref{LCH}, we have that 
\begin{align*}
\mathrm{sc}^{(0,1)}_{s}
\left(
\mathfrak{P}
\right)&=\mathrm{sc}^{(0,1)}_{s}
\left(
\mathfrak{Q}
\right),
&
\mathrm{tg}^{(0,1)}_{s}\left(
\mathfrak{P}
\right)&=\mathrm{tg}^{(0,1)}_{s}
\left(\mathfrak{Q}
\right).
\end{align*}
Hence, by Proposition~\ref{PPthCatAlg}, we have that
\begin{align*}
\mathrm{sc}_{s}^{0\mathbf{Pth}_{\boldsymbol{\mathcal{A}}}}
\left(
\mathfrak{P}
\right)
=
\mathrm{ip}^{(1,0)\sharp}_{s}
\left(
\mathrm{sc}^{(0,1)}_{s}
\left(
\mathfrak{P}
\right)\right)
=
\mathrm{ip}^{(1,0)\sharp}_{s}
\left(
\mathrm{sc}^{(0,1)}_{s}
\left(
\mathfrak{Q}
\right)\right)
=\mathrm{sc}_{s}^{0\mathbf{Pth}_{\boldsymbol{\mathcal{A}}}}
\left(
\mathfrak{Q}
\right),
\\
\mathrm{tg}_{s}^{0\mathbf{Pth}_{\boldsymbol{\mathcal{A}}}}
\left(
\mathfrak{P}
\right)
=
\mathrm{ip}^{(1,0)\sharp}_{s}
\left(
\mathrm{tg}^{(0,1)}_{s}
\left(
\mathfrak{P}
\right)\right)
=
\mathrm{ip}^{(1,0)\sharp}_{s}
\left(
\mathrm{tg}^{(0,1)}_{s}
\left(
\mathfrak{Q}
\right)\right)
=\mathrm{tg}_{s}^{0\mathbf{Pth}_{\boldsymbol{\mathcal{A}}}}
\left(
\mathfrak{Q}
\right).
\end{align*}

Since $\mathrm{sc}_{s}^{0\mathbf{Pth}_{\boldsymbol{\mathcal{A}}}}(\mathfrak{P})=\mathrm{sc}_{s}^{0\mathbf{Pth}_{\boldsymbol{\mathcal{A}}}}(\mathfrak{Q})$ and $\mathrm{tg}_{s}^{0\mathbf{Pth}_{\boldsymbol{\mathcal{A}}}}(\mathfrak{P})=\mathrm{tg}_{s}^{0\mathbf{Pth}_{\boldsymbol{\mathcal{A}}}}(\mathfrak{Q})$, it follows that they have, respectively, the same image under the Curry-Howard mapping. 
 
Therefore, we are only left with the compatibility of $\mathrm{Ker}(\mathrm{CH}^{(1)})$ with the $0$-composition operator. 

Let $s$ be a sort in $S$ and let $(\mathfrak{P},\mathfrak{Q})$ and $(\mathfrak{P}',\mathfrak{Q}')$ be two pairs of paths in $\mathrm{Ker}(\mathrm{CH}^{(1)})_{s}$. Let us assume that the $0$-composition 
$\mathfrak{P}'
\circ^{0\mathbf{Pth}_{\boldsymbol{\mathcal{A}}}}_{s}
\mathfrak{P}$ 
is defined. Then, by Lemma~\ref{LCH}, the composition 
$\mathfrak{Q}'
\circ^{0\mathbf{Pth}_{\boldsymbol{\mathcal{A}}}}_{s}
\mathfrak{Q}$ 
is also defined. We want to prove that their respective $0$-compositions, i.e., $\mathfrak{P}'\circ_{s}^{0\mathbf{Pth}_{\boldsymbol{\mathcal{A}}}}\mathfrak{P}$ and $\mathfrak{Q}'\circ_{s}^{0\mathbf{Pth}_{\boldsymbol{\mathcal{A}}}}\mathfrak{Q}$ have the same value for the Curry-Howard mapping.

We first study the case in which one of the paths is a $(1,0)$-identity path.

\begin{claim}\label{CCompIp} Let $s$ be a sort in $S$ and let $(\mathfrak{P},\mathfrak{Q}), (\mathfrak{P}',\mathfrak{Q}')$ be two pairs of paths in $\mathrm{Ker}(\mathrm{CH}^{(1)})_{s}$. If $\mathfrak{P},\mathfrak{P}',\mathfrak{Q}$ or $\mathfrak{Q}'$ is a $(1,0)$-identity path, then, when defined, 
$$
\left(
\mathfrak{P}'
\circ^{0\mathbf{Pth}_{\boldsymbol{\mathcal{A}}}}_{s}
\mathfrak{P}
,
\mathfrak{Q}'
\circ^{0\mathbf{Pth}_{\boldsymbol{\mathcal{A}}}}_{s}
\mathfrak{Q}
\right)
\in\mathrm{Ker}(\mathrm{CH}^{(1)})_{s}.$$
\end{claim}

Assume that $\mathfrak{P}$ is a $(1,0)$-identity path. 

Since $\mathfrak{P}'
\circ^{0\mathbf{Pth}_{\boldsymbol{\mathcal{A}}}}_{s}
\mathfrak{P}$ is defined, we have that $\mathrm{sc}^{(0,1)}_{s}(\mathfrak{P}')=\mathrm{tg}^{(0,1)}_{s}(\mathfrak{P})$. Hence $\mathfrak{P}$ is the $(1,0)$-identity path on the $(0,1)$-source of $\mathfrak{P}'$, i.e., $\mathfrak{P}=\mathrm{ip}^{(1,0)\sharp}_{s}(\mathrm{sc}^{(0,1)}_{s}(\mathfrak{P}'))$. Hence, the $0$-composition $\mathfrak{P}'
\circ^{0\mathbf{Pth}_{\boldsymbol{\mathcal{A}}}}_{s}
\mathfrak{P}$  reduces to $\mathfrak{P}'$.

On the other hand, since $(\mathfrak{P},\mathfrak{Q})$ is a pair in $\mathrm{Ker}(\mathrm{CH}^{(1)})_{s}$, we have, by Corollary~\ref{CCHUZId}, that $\mathfrak{P}=\mathfrak{Q}$. Hence, $\mathfrak{Q}=\mathrm{ip}^{(1,0)\sharp}_{s}(\mathrm{sc}^{(0,1)}_{s}(\mathfrak{P}'))$. Moreover, since $(\mathfrak{P}',\mathfrak{Q}')$ is a pair in $\mathrm{Ker}(\mathrm{CH}^{(1)})_{s}$, it follows that, by Lemma~\ref{LCH}, $\mathrm{sc}^{(0,1)}_{s}(\mathfrak{P}')=\mathrm{sc}^{(0,1)}_{s}(\mathfrak{Q}')$. Hence, the $0$-composition $\mathfrak{Q}'
\circ^{0\mathbf{Pth}_{\boldsymbol{\mathcal{A}}}}_{s}
\mathfrak{Q}$  reduces to $\mathfrak{Q}'$.

Then we have that
$$
\left(
\mathfrak{P}'
\circ^{0\mathbf{Pth}_{\boldsymbol{\mathcal{A}}}}_{s}
\mathfrak{P}
,
\mathfrak{Q}'
\circ^{0\mathbf{Pth}_{\boldsymbol{\mathcal{A}}}}_{s}
\mathfrak{Q}
\right)
\in\mathrm{Ker}(\mathrm{CH}^{(1)})_{s}.$$

The remaining cases are proved similarly. This finishes the proof of Claim~\ref{CCompIp}.

We are now in position to prove the general situation.

\begin{claim}\label{CCompProd} Let $s$ be a sort in $S$ and let $(\mathfrak{P},\mathfrak{Q}), (\mathfrak{P}',\mathfrak{Q}')$ be two pairs of paths in $\mathrm{Ker}(\mathrm{CH}^{(1)})_{s}$. When defined, 
$$
\left(
\mathfrak{P}'
\circ^{0\mathbf{Pth}_{\boldsymbol{\mathcal{A}}}}_{s}
\mathfrak{P}
,
\mathfrak{Q}'
\circ^{0\mathbf{Pth}_{\boldsymbol{\mathcal{A}}}}_{s}
\mathfrak{Q}
\right)
\in\mathrm{Ker}(\mathrm{CH}^{(1)})_{s}.$$
\end{claim}

The proof is done by Artinian induction on $(\coprod\mathrm{Pth}_{\boldsymbol{\mathcal{A}}},\leq_{\mathbf{Pth}_{\boldsymbol{\mathcal{A}}}})$.

\textsf{Base step of the Artinian induction}.

Let $(\mathfrak{P}'\circ^{0\mathbf{Pth}_{\boldsymbol{\mathcal{A}}}}_{s}\mathfrak{P},s)$ be a minimal element in $(\coprod\mathrm{Pth}_{\boldsymbol{\mathcal{A}}}, \leq_{\mathbf{Pth}_{\boldsymbol{\mathcal{A}}}})$. Then, by Proposition~\ref{PMinimal}, the path $\mathfrak{P}'\circ^{0\mathbf{Pth}_{\boldsymbol{\mathcal{A}}}}_{s}\mathfrak{P}$ is either a $(1,0)$-identity path on a simple term, or an echelon. In any case, either $\mathfrak{P}$ or $\mathfrak{P}'$ must be a $(1,0)$-identity path. The statement follows by Claim~\ref{CCompIp}.

\textsf{Inductive step of the Artinian induction}.

Let $(\mathfrak{P}'\circ^{0\mathbf{Pth}_{\boldsymbol{\mathcal{A}}}}_{s}\mathfrak{P},s)$ be a non-minimal element in $(\coprod\mathrm{Pth}_{\boldsymbol{\mathcal{A}}}, \leq_{\mathbf{Pth}_{\boldsymbol{\mathcal{A}}}})$. Let us suppose that, for every sort $t\in S$ and every path $\mathfrak{P}'''\circ^{0\mathbf{Pth}_{\boldsymbol{\mathcal{A}}}}_{t}\mathfrak{P}''\in\mathrm{Pth}_{\boldsymbol{\mathcal{A}},t}$, if $(\mathfrak{P}'''\circ^{0\mathbf{Pth}_{\boldsymbol{\mathcal{A}}}}_{t}\mathfrak{P}'',t)<_{\mathbf{Pth}_{\boldsymbol{\mathcal{A}}}}(\mathfrak{P}'\circ^{0\mathbf{Pth}_{\boldsymbol{\mathcal{A}}}}_{s}\mathfrak{P},s)$, then the statement holds for $\mathfrak{P}'''\circ^{0\mathbf{Pth}_{\boldsymbol{\mathcal{A}}}}_{t}\mathfrak{P}''$, i.e., for every pair of paths $\mathfrak{Q}''$ and $\mathfrak{Q}'''$ in $\mathrm{Pth}_{\boldsymbol{\mathcal{A}},t}$, if $(\mathfrak{P}'',\mathfrak{Q}'')\in\mathrm{Ker}(\mathrm{CH}^{(1)})_{t}$ and $(\mathfrak{P}''',\mathfrak{Q}''')\in\mathrm{Ker}(\mathrm{CH}^{(1)})_{t}$,  then 
$$
\left(
\mathfrak{P}'''\circ^{0\mathbf{Pth}_{\boldsymbol{\mathcal{A}}}}_{t}\mathfrak{P}'',\mathfrak{Q}'''\circ^{0\mathbf{Pth}_{\boldsymbol{\mathcal{A}}}}_{t}\mathfrak{Q}''
\right)\in \mathrm{Ker}(\mathrm{CH}^{(1)})_{t}.$$

Since $(\mathfrak{P}'\circ^{0\mathbf{Pth}_{\boldsymbol{\mathcal{A}}}}_{s}\mathfrak{P},s)$ is a non-minimal element in $(\coprod\mathrm{Pth}_{\boldsymbol{\mathcal{A}}}, \leq_{\mathbf{Pth}_{\boldsymbol{\mathcal{A}}}})$ and, by Claim~\ref{CCompIp}, we can assume that neither $\mathfrak{P}'$ nor $\mathfrak{P}$ are $(1,0)$-identity paths, we have, by Lemma~\ref{LOrdI}, that $\mathfrak{P}'\circ^{0\mathbf{Pth}_{\boldsymbol{\mathcal{A}}}}_{s}\mathfrak{P}$ is either~(1) a path of length strictly greater than one containing at least one echelon or~(2) an echelonless path. 

If~(1), then let $i\in \bb{\mathfrak{P}'\circ^{0\mathbf{Pth}_{\boldsymbol{\mathcal{A}}}}_{s}\mathfrak{P}}$ be the first index for which the one-step subpath $(\mathfrak{P}'\circ^{0\mathbf{Pth}_{\boldsymbol{\mathcal{A}}}}_{s}\mathfrak{P})^{i,i}$ of $\mathfrak{P}'\circ^{0\mathbf{Pth}_{\boldsymbol{\mathcal{A}}}}_{s}\mathfrak{P}$ is an echelon. We distinguish the cases~(1.1) $i=0$ and~(1.2) $i>0$.

If~(1.1), i.e., $i=0$, since we are assuming that $\mathfrak{P}$ is not a $(1,0)$-identity path, we have that $\mathfrak{P}$ has an echelon in its first step. Then it could be the case that either (1.1.1)~$\mathfrak{P}$ is an echelon or (1.1.2)~$\mathfrak{P}$ is a path of length strictly greater than one containing an echelon on its first step.

If~(1.1.1) then, since $(\mathfrak{P},\mathfrak{Q})\in\mathrm{Ker}(\mathrm{CH}^{(1)})_{s}$, we have, by Proposition~\ref{PCHEch}, that $\mathfrak{P}=\mathfrak{Q}$. Thus, the paths $\mathfrak{P}'\circ^{0\mathbf{Pth}_{\boldsymbol{\mathcal{A}}}}_{s}\mathfrak{P}$ and $\mathfrak{Q}'\circ^{0\mathbf{Pth}_{\boldsymbol{\mathcal{A}}}}_{s}\mathfrak{Q}$ are both paths of length strictly greater than one containing an echelon in its first step, this echelon being, respectively, $\mathfrak{P}$ and $\mathfrak{Q}$. Then in virtue of Proposition~\ref{PCHHomCat} we have that 
\begin{align*}
\mathrm{CH}^{(1)}_{s}
\left(
\mathfrak{P}'\circ^{0\mathbf{Pth}_{\boldsymbol{\mathcal{A}}}}_{s}\mathfrak{P}
\right)
&=
\mathrm{CH}^{(1)}_{s}
\left(\mathfrak{P}'
\right)
\circ^{0
\mathbf{T}_{\Sigma^{\boldsymbol{\mathcal{A}}}}(X)
}_{s}
\mathrm{CH}^{(1)}_{s}
\left(
\mathfrak{P}
\right)
\\&=
\mathrm{CH}^{(1)}_{s}
\left(\mathfrak{Q}'
\right)
\circ^{0
\mathbf{T}_{\Sigma^{\boldsymbol{\mathcal{A}}}}(X)
}_{s}
\mathrm{CH}^{(1)}_{s}
\left(
\mathfrak{Q}
\right)
\\&=
\mathrm{CH}^{(1)}_{s}
\left(
\mathfrak{Q}'
\circ^{0\mathbf{Pth}_{\boldsymbol{\mathcal{A}}}}_{s}\mathfrak{Q}
\right).
\end{align*}

If~(1.1.2) then the value of the Curry-Howard mapping at $\mathfrak{P}$ is given by
$$
\mathrm{CH}^{(1)}_{s}\left(
\mathfrak{P}
\right)
=
\mathrm{CH}^{(1)}_{s}
\left(
\mathfrak{P}^{1,\bb{\mathfrak{P}}-1}
\right)
\circ^{0
\mathbf{T}_{\Sigma^{\boldsymbol{\mathcal{A}}}}(X)
}_{s}
\mathrm{CH}^{(1)}_{s}\left(
\mathfrak{P}^{0,0}
\right).
$$

Since $(\mathfrak{P},\mathfrak{Q})\in \mathrm{Ker}(\mathrm{CH}^{(1)})_{s}$ and $\mathrm{CH}^{(1)}_{s}(\mathfrak{P})\in\eta^{(1,\mathcal{A})}[\mathcal{A}]^{\mathrm{int}}_{s}$, we have, by Lemma~\ref{LCHEchInt}, that $\mathfrak{Q}$ is a path of length strictly greater than one containing an echelon on its first step.  The same applies for the the path $\mathfrak{Q}'\circ_{s}^{0\mathbf{Pth}_{\boldsymbol{\mathcal{A}}}}\mathfrak{Q}$.

Then, the value of the Curry-Howard mapping at $\mathfrak{Q}$ is given by
$$
\mathrm{CH}^{(1)}_{s}\left(
\mathfrak{Q}
\right)
=
\mathrm{CH}^{(1)}_{s}\left(
\mathfrak{Q}^{1,\bb{\mathfrak{Q}}-1}
\right)
\circ^{0
\mathbf{T}_{\Sigma^{\boldsymbol{\mathcal{A}}}}(X)
}_{s}
\mathrm{CH}^{(1)}_{s}\left(
\mathfrak{Q}^{0,0}
\right).
$$

Since $(\mathfrak{P},\mathfrak{Q})\in \mathrm{Ker}(\mathrm{CH}^{(1)})_{s}$, it follows that  $(\mathfrak{P}^{1,\bb{\mathfrak{P}}-1},\mathfrak{Q}^{1,\bb{\mathfrak{Q}}-1})$ and $(\mathfrak{P}^{0,0},\mathfrak{Q}^{0,0})$ are pairs in $\mathrm{Ker}(\mathrm{CH}^{(1)})_{s}$. Moreover, taking into account that 
\begin{align*}
\mathrm{tg}^{(0,1)}_{s}\left(\mathfrak{Q}
\right)&=
\mathrm{tg}^{(0,1)}_{s}\left(
\mathfrak{Q}^{1,\bb{\mathfrak{Q}}-1}
\right),
&
\mathrm{tg}^{(0,1)}_{s}\left(
\mathfrak{P}
\right)&=
\mathrm{tg}^{(0,1)}_{s}
\left(
\mathfrak{P}^{1,\bb{\mathfrak{P}}-1}
\right),
\end{align*}
the following $0$-compositions are well-defined
\begin{align*}
\mathfrak{P}'\circ_{s}^{0\mathbf{Pth}_{\boldsymbol{\mathcal{A}}}}
 \mathfrak{P}^{1,\bb{\mathfrak{P}}-1},
&&
\mathfrak{Q}'\circ_{s}^{0\mathbf{Pth}_{\boldsymbol{\mathcal{A}}}}
 \mathfrak{Q}^{1,\bb{\mathfrak{Q}}-1}.
\end{align*}

Since $(\mathfrak{P}'\circ_{s}^{0\mathbf{Pth}_{\boldsymbol{\mathcal{A}}}}
 \mathfrak{P}^{1,\bb{\mathfrak{P}}-1},s)\prec_{\mathbf{Pth}_{\boldsymbol{\mathcal{A}}}}(\mathfrak{P}'\circ_{s}^{0\mathbf{Pth}_{\boldsymbol{\mathcal{A}}}}\mathfrak{P},s)$, and the pairs $(\mathfrak{P}^{1,\bb{\mathfrak{P}}-1},\mathfrak{Q}^{1,\bb{\mathfrak{Q}}-1})$ and $(\mathfrak{P}',\mathfrak{Q}')$ are in $\mathrm{Ker}(\mathrm{CH}^{(1)})_{s}$, we have, by induction, that
$$
\left(\mathfrak{P}'\circ_{s}^{0\mathbf{Pth}_{\boldsymbol{\mathcal{A}}}}
 \mathfrak{P}^{1,\bb{\mathfrak{P}}-1},
\mathfrak{Q}'\circ_{s}^{0\mathbf{Pth}_{\boldsymbol{\mathcal{A}}}}
 \mathfrak{Q}^{1,\bb{\mathfrak{Q}}-1}
 \right)
 \in\mathrm{Ker}(\mathrm{CH}^{(1)})_{s}. 
$$

So, considering the foregoing, we can affirm that
\begin{flushleft}
$\mathrm{CH}^{(1)}_{s}\left(
\mathfrak{P}'\circ_{s}^{0\mathbf{Pth}_{\boldsymbol{\mathcal{A}}}}
\mathfrak{P}
\right)$
\begin{align*}
\qquad
&=
\mathrm{CH}^{(1)}_{s}
\left(\left(
\mathfrak{P}'\circ_{s}^{0\mathbf{Pth}_{\boldsymbol{\mathcal{A}}}}
\mathfrak{P}
\right)^{1,\bb{\mathfrak{P}'\circ_{s}^{0\mathbf{Pth}_{\boldsymbol{\mathcal{A}}}}
\mathfrak{P}}-1}
\right)
\circ_{s}^{0\mathbf{T}_{\Sigma^{\boldsymbol{\mathcal{A}}}}(X)}
\mathrm{CH}^{(1)}_{s}
\left(\left(
\mathfrak{P}'\circ_{s}^{0\mathbf{Pth}_{\boldsymbol{\mathcal{A}}}}
\mathfrak{P}
\right)^{0,0}
\right)
\\&=
\mathrm{CH}^{(1)}_{s}
\left(
\mathfrak{P}'\circ_{s}^{0\mathbf{Pth}_{\boldsymbol{\mathcal{A}}}}
\mathfrak{P}^{1,\bb{
\mathfrak{P}}-1}
\right)
\circ_{s}^{0\mathbf{T}_{\Sigma^{\boldsymbol{\mathcal{A}}}}(X)}
\mathrm{CH}^{(1)}_{s}
\left(
\mathfrak{P}^{0,0}
\right)
\\&=
\mathrm{CH}^{(1)}_{s}
\left(
\mathfrak{Q}'\circ_{s}^{0\mathbf{Pth}_{\boldsymbol{\mathcal{A}}}}
\mathfrak{Q}^{1,\bb{
\mathfrak{Q}}-1}
\right)
\circ_{s}^{0\mathbf{T}_{\Sigma^{\boldsymbol{\mathcal{A}}}}(X)}
\mathrm{CH}^{(1)}_{s}
\left(
\mathfrak{Q}^{0,0}
\right)
\\&=
\mathrm{CH}^{(1)}_{s}
\left(\left(
\mathfrak{Q}'\circ_{s}^{0\mathbf{Pth}_{\boldsymbol{\mathcal{A}}}}
\mathfrak{Q}
\right)^{1,\bb{\mathfrak{Q}'\circ_{s}^{0\mathbf{Pth}_{\boldsymbol{\mathcal{A}}}}
\mathfrak{Q}}-1}
\right)
\circ_{s}^{0\mathbf{T}_{\Sigma^{\boldsymbol{\mathcal{A}}}}(X)}
\mathrm{CH}^{(1)}_{s}
\left(\left(
\mathfrak{Q}'\circ_{s}^{0\mathbf{Pth}_{\boldsymbol{\mathcal{A}}}}
\mathfrak{Q}
\right)^{0,0}
\right)
\\&=
\mathrm{CH}^{(1)}_{s}
\left(
\mathfrak{Q}'\circ_{s}^{0\mathbf{Pth}_{\boldsymbol{\mathcal{A}}}}
\mathfrak{Q}
\right).
\end{align*}
\end{flushleft}

This finishes the case $i=0$.

For case (1.2), i.e., if $i\neq 0$, since $\bb{\mathfrak{P}'\circ^{0\mathbf{Pth}_{\boldsymbol{\mathcal{A}}}}_{s}\mathfrak{P}}=\bb{\mathfrak{P}'}+\bb{\mathfrak{P}}$, then either~(1.2.1) $i\in\bb{\mathfrak{P}}$ or~(1.2.2) $i\in[\bb{\mathfrak{P}},\bb{\mathfrak{P}'\circ^{0\mathbf{Pth}_{\boldsymbol{\mathcal{A}}}}_{s}\mathfrak{P}}-1]$. 

If~(1.2.1), i.e., $i\neq 0$ and $i\in\bb{\mathfrak{P}}$, then $\mathfrak{P}$ is a path of length strictly greater than one containing an echelon in a step different from the initial one.

Then the value of the Curry-Howard mapping at $\mathfrak{P}$ is given by
$$
\mathrm{CH}^{(1)}_{s}\left(
\mathfrak{P}
\right)=
\mathrm{CH}^{(1)}_{s}
\left(
\mathfrak{P}^{i,\bb{\mathfrak{P}}-1}
\right)
\circ_{s}^{0\mathbf{T}_{\Sigma^{\boldsymbol{\mathcal{A}}}}(X)}
\mathrm{CH}^{(1)}_{s}
\left(
\mathfrak{P}^{0,i-1}
\right).
$$

Since $(\mathfrak{P},\mathfrak{Q})\in \mathrm{Ker}(\mathrm{CH}^{(1)})_{s}$, and $\mathrm{CH}^{(1)}_{s}(\mathfrak{P})\in\eta^{(1,\mathcal{A})}[\mathcal{A}]^{\neg{\mathrm{int}}}_{s}$, we have, by Lemma~\ref{LCHEchNInt}, that $\mathfrak{Q}$ is a path of length strictly greater than one containing an echelon on a step different from the initial one.  The same applies for the the path $\mathfrak{Q}'\circ_{s}^{0\mathbf{Pth}_{\boldsymbol{\mathcal{A}}}}\mathfrak{Q}$.

Let $j\in\bb{\mathfrak{Q}}$ be the first index for which $\mathfrak{Q}^{j,j}$ is an echelon. By the previous discussion $j\neq 0$. Then the value of the Curry-Howard mapping at $\mathfrak{Q}$ is given by
$$
\mathrm{CH}^{(1)}_{s}\left(
\mathfrak{Q}
\right)
=
\mathrm{CH}^{(1)}_{s}
\left(
\mathfrak{Q}^{j,\bb{\mathfrak{Q}}-1}
\right)
\circ^{0
\mathbf{T}_{\Sigma^{\boldsymbol{\mathcal{A}}}}(X)
}_{s}
\mathrm{CH}^{(1)}_{s}
\left(
\mathfrak{Q}^{0,j-1}
\right).
$$

Since $(\mathfrak{P},\mathfrak{Q})\in \mathrm{Ker}(\mathrm{CH}^{(1)})_{s}$, we have that  $(\mathfrak{P}^{i,\bb{\mathfrak{P}}-1},\mathfrak{Q}^{j,\bb{\mathfrak{Q}}-1})$ and $(\mathfrak{P}^{0,i-1},\mathfrak{Q}^{0,j-1})$ are pairs in $\mathrm{Ker}(\mathrm{CH}^{(1)})_{s}$. Moreover, taking into account that 
\begin{align*}
\mathrm{tg}^{(0,1)}_{s}\left(
\mathfrak{Q}
\right)&=
\mathrm{tg}^{(0,1)}_{s}
\left(
\mathfrak{Q}^{j,\bb{\mathfrak{Q}}-1}
\right)
&
\mathrm{tg}^{(0,1)}_{s}
\left(
\mathfrak{P}
\right)&=
\mathrm{tg}^{(0,1)}_{s}
\left(
\mathfrak{P}^{i,\bb{\mathfrak{P}}-1}
\right),
\end{align*}
the following $0$-compositions are well-defined
\begin{align*}
\mathfrak{P}'\circ_{s}^{0\mathbf{Pth}_{\boldsymbol{\mathcal{A}}}}
 \mathfrak{P}^{i,\bb{\mathfrak{P}}-1},
&&
\mathfrak{Q}'\circ_{s}^{0\mathbf{Pth}_{\boldsymbol{\mathcal{A}}}}
 \mathfrak{Q}^{j,\bb{\mathfrak{Q}}-1}.
\end{align*}

Since $(\mathfrak{P}'\circ_{s}^{0\mathbf{Pth}_{\boldsymbol{\mathcal{A}}}}
 \mathfrak{P}^{i,\bb{\mathfrak{P}}-1},s)\prec_{\mathbf{Pth}_{\boldsymbol{\mathcal{A}}}}(\mathfrak{P}'\circ_{s}^{0\mathbf{Pth}_{\boldsymbol{\mathcal{A}}}}\mathfrak{P},s)$, and the pairs $(\mathfrak{P}^{i,\bb{\mathfrak{P}}-1},\mathfrak{Q}^{j,\bb{\mathfrak{Q}}-1})$ and $(\mathfrak{P}',\mathfrak{Q}')$ are in $\mathrm{Ker}(\mathrm{CH}^{(1)})_{s}$, we have, by induction,  that
$$
\left(\mathfrak{P}'\circ_{s}^{0\mathbf{Pth}_{\boldsymbol{\mathcal{A}}}}
 \mathfrak{P}^{i,\bb{\mathfrak{P}}-1},
\mathfrak{Q}'\circ_{s}^{0\mathbf{Pth}_{\boldsymbol{\mathcal{A}}}}
 \mathfrak{Q}^{j,\bb{\mathfrak{Q}}-1}
 \right)
 \in\mathrm{Ker}(\mathrm{CH}^{(1)})_{s}. 
$$

So, considering the foregoing, we can affirm that
\begin{flushleft}
$\mathrm{CH}^{(1)}_{s}
\left(
\mathfrak{P}'\circ_{s}^{0\mathbf{Pth}_{\boldsymbol{\mathcal{A}}}}
\mathfrak{P}
\right)$
\begin{align*}
\qquad
&=
\mathrm{CH}^{(1)}_{s}
\left(\left(
\mathfrak{P}'\circ_{s}^{0\mathbf{Pth}_{\boldsymbol{\mathcal{A}}}}
\mathfrak{P}
\right)^{i,\bb{\mathfrak{P}'\circ_{s}^{0\mathbf{Pth}_{\boldsymbol{\mathcal{A}}}}
\mathfrak{P}}-1}
\right)
\circ_{s}^{0\mathbf{T}_{\Sigma^{\boldsymbol{\mathcal{A}}}}(X)}
\mathrm{CH}^{(1)}_{s}
\left(\left(\mathfrak{P}'\circ_{s}^{0\mathbf{Pth}_{\boldsymbol{\mathcal{A}}}}
\mathfrak{P}
\right)^{0,i-1}
\right)
\\&=
\mathrm{CH}^{(1)}_{s}
\left(
\mathfrak{P}'\circ_{s}^{0\mathbf{Pth}_{\boldsymbol{\mathcal{A}}}}
\mathfrak{P}^{i,\bb{
\mathfrak{P}}-1}
\right)
\circ_{s}^{0\mathbf{T}_{\Sigma^{\boldsymbol{\mathcal{A}}}}(X)}
\mathrm{CH}^{(1)}_{s}
\left(
\mathfrak{P}^{0,i-1}
\right)
\\&=
\mathrm{CH}^{(1)}_{s}
\left(
\mathfrak{Q}'\circ_{s}^{0\mathbf{Pth}_{\boldsymbol{\mathcal{A}}}}
\mathfrak{Q}^{j,\bb{
\mathfrak{Q}}-1}
\right)
\circ_{s}^{0\mathbf{T}_{\Sigma^{\boldsymbol{\mathcal{A}}}}(X)}
\mathrm{CH}^{(1)}_{s}
\left(
\mathfrak{Q}^{0,j-1}
\right)
\\&=
\mathrm{CH}^{(1)}_{s}
\left(\left(
\mathfrak{Q}'\circ_{s}^{0\mathbf{Pth}_{\boldsymbol{\mathcal{A}}}}
\mathfrak{Q}
\right)^{j,\bb{\mathfrak{Q}'\circ_{s}^{0\mathbf{Pth}_{\boldsymbol{\mathcal{A}}}}
\mathfrak{Q}}-1}
\right)
\circ_{s}^{0\mathbf{T}_{\Sigma^{\boldsymbol{\mathcal{A}}}}(X)}
\mathrm{CH}^{(1)}_{s}
\left(\left(\mathfrak{Q}'\circ_{s}^{0\mathbf{Pth}_{\boldsymbol{\mathcal{A}}}}
\mathfrak{Q}
\right)^{0,j-1}
\right)
\\&=
\mathrm{CH}^{(1)}_{s}\left(
\mathfrak{Q}'\circ_{s}^{0\mathbf{Pth}_{\boldsymbol{\mathcal{A}}}}
\mathfrak{Q}
\right).
\end{align*}
\end{flushleft}

The case $i\in\bb{\mathfrak{P}}$ follows.

If~(1.2.2), i.e. $i\neq 0$ and $i\in [\bb{\mathfrak{P}}, \bb{
\mathfrak{P}'\circ^{0\mathbf{Pth}_{\boldsymbol{\mathcal{A}}}}_{s}\mathfrak{P}
}-1]$, then  $\mathfrak{P}'$ is a non $(1,0)$-identity path containing an echelon, whilst $\mathfrak{P}$ is an echelonless path.

We will distinguish three cases according to whether~(1.2.2.1) $\mathfrak{P}'$ is an echelon;~(1.2.2.2) $\mathfrak{P}'$ is a path of length strictly greater than one containing an echelon on its first step or~(1.2.2.3) $\mathfrak{P}'$ is a path of length strictly greater than one containing an echelon on a step different from zero. These cases can be proved using a similar argument to those three cases presented above. We leave the details for the interested reader.

If~(2), i.e., if $\mathfrak{P}'\circ^{0\mathbf{Pth}_{\boldsymbol{\mathcal{A}}}}_{s}\mathfrak{P}$ is an echelonless path then, regarding the paths $\mathfrak{P}'$ and $\mathfrak{P}$, we have that
\begin{itemize}
\item[(i)] $\mathfrak{P}$ is an echelonless path.
\item[(ii)] $\mathfrak{P}'$ is an echelonless path.
\end{itemize}

Since~(i), we have by Lemma~\ref{LPthHeadCt}, that there exists a unique word $\mathbf{s}\in S^{\star}-\{\lambda\}$ and a unique  operation symbol $\sigma\in \Sigma_{\mathbf{s},s}$ associated to $\mathfrak{P}$. Let $(\mathfrak{P}_{j})_{j\in\bb{\mathbf{s}}}$ be the family of paths we can extract from $\mathfrak{P}$ in virtue of Lemma~\ref{LPthExtract}. Then, according to Definition~\ref{DCH}, we have that the value of the Curry-Howard mapping at $\mathfrak{P}$ is given by 
$$
\mathrm{CH}^{(1)}_{s}\left(
\mathfrak{P}
\right)
=
\sigma^{\mathbf{T}_{\Sigma^{\boldsymbol{\mathcal{A}}}}(X)}
\left(\left(
\mathrm{CH}^{(1)}_{s_{j}}\left(
\mathfrak{P}_{j}
\right)\right)_{j\in\bb{\mathbf{s}}}
\right).
$$

Since $(\mathfrak{P},\mathfrak{Q})\in\mathrm{Ker}(\mathrm{CH}^{(1)})_{s}$ and $\mathrm{CH}^{(1)}_{s}(\mathfrak{P})\in\mathcal{T}(\sigma,\mathrm{T}_{\Sigma^{\boldsymbol{\mathcal{A}}}}(X))_{1}$ we have, by Lemma~\ref{LCHNEch} that
\begin{itemize}
\item[(iii)] $\mathfrak{Q}$ is an echelonless path associated to the operation symbol $\sigma\in\Sigma_{\mathbf{s},s}$.
\end{itemize}

Since (iii), let $(\mathfrak{Q}_{j})_{j\in\bb{\mathbf{s}}}$ be the family of paths we can extract from $\mathfrak{Q}$ in virtue of Lemma~\ref{LPthExtract}. Then, according to Definition~\ref{DCH}, we have that the value of the Curry-Howard mapping at $\mathfrak{Q}$ is given by
$$
\mathrm{CH}^{(1)}_{s}\left(
\mathfrak{Q}
\right)
=
\sigma^{\mathbf{T}_{\Sigma^{\boldsymbol{\mathcal{A}}}}(X)}
\left(\left(
\mathrm{CH}^{(1)}_{s_{j}}\left(
\mathfrak{Q}_{j}
\right)\right)_{j\in\bb{\mathbf{s}}}
\right).
$$

Since $(\mathfrak{P},\mathfrak{Q})$ is a pair in $\mathrm{Ker}(\mathrm{CH}^{(1)})_{s}$, then we have that, for every $j\in\bb{\mathbf{s}}$, it happens that 
$$
\left(
\mathfrak{P}_{j},
\mathfrak{Q}_{j}
\right)
\in\mathrm{Ker}\left(
\mathrm{CH}^{(1)}
\right)_{s_{j}}.
$$

Since~(ii), we have by Lemma~\ref{LPthHeadCt}, that for the unique word $\mathbf{s}\in S^{\star}-\{\lambda\}$ and the unique  operation symbol $\sigma\in \Sigma_{\mathbf{s},s}$,  $\sigma$ is the operation symbol associated to $\mathfrak{P}'$.

Note that the operation symbol $\sigma$ is the same as in case (i), since $\mathfrak{P}'\circ_{s}^{0\mathbf{Pth}_{\boldsymbol{\mathcal{A}}}}\mathfrak{P}$ is an echelonless path by hypothesis and, thus, head-constant by Lemma~\ref{LPthHeadCt}.

Let $(\mathfrak{P}'_{j})_{j\in\bb{\mathbf{s}}}$ be the family of paths we can extract from $\mathfrak{P}'$ in virtue of Lemma~\ref{LPthExtract}. Then, according to Definition~\ref{DCH}, we have that the value of the Curry-Howard mapping at $\mathfrak{P}'$ is given by 
$$
\mathrm{CH}^{(1)}_{s}\left(
\mathfrak{P}'
\right)
=
\sigma^{\mathbf{T}_{\Sigma^{\boldsymbol{\mathcal{A}}}}(X)}
\left(\left(
\mathrm{CH}^{(1)}_{s_{j}}\left(
\mathfrak{P}'_{j}
\right)\right)_{j\in\bb{\mathbf{s}}}
\right).
$$

Since $(\mathfrak{P}',\mathfrak{Q}')\in\mathrm{Ker}(\mathrm{CH}^{(1)})_{s}$ and $\mathrm{CH}^{(1)}_{s}(\mathfrak{P}')\in\mathcal{T}(\sigma,\mathrm{T}_{\Sigma^{\boldsymbol{\mathcal{A}}}}(X))_{1}$ we have, by Lemma~\ref{LCHNEch} that
\begin{itemize}
\item[(iv)] $\mathfrak{Q}'$ is an echelonless path associated to the operation symbol $\sigma\in\Sigma_{\mathbf{s},s}$.
\end{itemize}

Since (iv), let $(\mathfrak{Q}'_{j})_{j\in\bb{\mathbf{s}}}$ be the family of paths we can extract from $\mathfrak{Q}'$ in virtue of Lemma~\ref{LPthExtract}. Then, according to Definition~\ref{DCH}, we have that the value of the Curry-Howard mapping at $\mathfrak{Q}'$ is given by
$$
\mathrm{CH}^{(1)}_{s}\left(
\mathfrak{Q}'
\right)
=
\sigma^{\mathbf{T}_{\Sigma^{\boldsymbol{\mathcal{A}}}}(X)}
\left(\left(
\mathrm{CH}^{(1)}_{s_{j}}\left(
\mathfrak{Q}'_{j}
\right)\right)_{j\in\bb{\mathbf{s}}}
\right).
$$

Since $(\mathfrak{P}',\mathfrak{Q}')$ is a pair in $\mathrm{Ker}(\mathrm{CH}^{(1)})_{s}$, then we have that, for every $j\in\bb{\mathbf{s}}$, it happens that 
$$
\left(
\mathfrak{P}'_{j},
\mathfrak{Q}'_{j}
\right)
\in\mathrm{Ker}\left(
\mathrm{CH}^{(1)}
\right)_{s_{j}}.
$$

From (iii) and (iv), we infer that $\mathfrak{Q}'\circ_{s}^{0\mathbf{Pth}_{\boldsymbol{\mathcal{A}}}}\mathfrak{Q}$ is an echelonless path. 

Let us consider $((\mathfrak{P}'\circ_{s}^{0\mathbf{Pth}_{\boldsymbol{\mathcal{A}}}}\mathfrak{P})_{j})_{j\in\bb{\mathbf{s}}}$ and $((\mathfrak{Q}'\circ_{s}^{0\mathbf{Pth}_{\boldsymbol{\mathcal{A}}}}\mathfrak{Q})_{j})_{j\in\bb{\mathbf{s}}}$ the family of paths we can extract, in virtue of Lemma~\ref{LPthExtract}, from $\mathfrak{P}'\circ_{s}^{0\mathbf{Pth}_{\boldsymbol{\mathcal{A}}}}\mathfrak{P}$ and $\mathfrak{Q}'\circ_{s}^{0\mathbf{Pth}_{\boldsymbol{\mathcal{A}}}}\mathfrak{Q}$, respectively. Let us note that, for every $j\in\bb{\mathbf{s}}$, it is the case that 
\begin{align*}
\left(
\mathfrak{P}'\circ_{s}^{0\mathbf{Pth}_{\boldsymbol{\mathcal{A}}}}\mathfrak{P}
\right)_{j}
&=\mathfrak{P}'_{j}\circ_{s_{j}}^{0\mathbf{Pth}_{\boldsymbol{\mathcal{A}}}}\mathfrak{P}_{j};
\\
\left(
\mathfrak{Q}'\circ_{s}^{0\mathbf{Pth}_{\boldsymbol{\mathcal{A}}}}\mathfrak{Q}
\right)_{j}
&=\mathfrak{Q}'_{j}\circ_{s_{j}}^{0\mathbf{Pth}_{\boldsymbol{\mathcal{A}}}}\mathfrak{Q}_{j}.
\end{align*}

For every $j\in\bb{\mathbf{s}}$, $(\mathfrak{P}'_{j}\circ_{s_{j}}^{0\mathbf{Pth}_{\boldsymbol{\mathcal{A}}}}\mathfrak{P}_{j},s_{j})
\prec_{\mathbf{Pth}_{\boldsymbol{\mathcal{A}}}}
(\mathfrak{P}'\circ_{s}^{0\mathbf{Pth}_{\boldsymbol{\mathcal{A}}}}\mathfrak{P},s)
$ and $(\mathfrak{P}'_{j},\mathfrak{Q}'_{j})$ and $(\mathfrak{P}_{j},\mathfrak{Q}_{j})$ are pairs in $\mathrm{Ker}(\mathrm{CH}^{(1)})_{s_{j}}$, hence, by induction, we have that, for every $j\in\bb{\mathbf{s}}$,
$$
\left(
\mathfrak{P}'_{j}
\circ_{s_{j}}^{0\mathbf{Pth}_{\boldsymbol{\mathcal{A}}}}
\mathfrak{P}_{j},
\mathfrak{Q}'_{j}
\circ_{s_{j}}^{0\mathbf{Pth}_{\boldsymbol{\mathcal{A}}}}
\mathfrak{Q}_{j}
\right)
\in\mathrm{Ker}
\left(
\mathrm{CH}^{(1)}
\right)_{s_{j}}.
$$

So, considering the foregoing, we can affirm that
\allowdisplaybreaks
\begin{align*}
\mathrm{CH}^{(1)}_{s}\left(
\mathfrak{P}'
\circ_{s}^{0\mathbf{Pth}_{\boldsymbol{\mathcal{A}}}}
\mathfrak{P}
\right)
&=
\sigma^{\mathrm{T}_{\Sigma^{\boldsymbol{\mathcal{A}}}}(X)}
\left(\left(
\mathrm{CH}^{(1)}_{s_{j}}\left(\left(
\mathfrak{P}'
\circ_{s}^{0\mathbf{Pth}_{\boldsymbol{\mathcal{A}}}}
\mathfrak{P}
\right)_{j}
\right)\right)_{j\in\bb{\mathbf{s}}}
\right)
\\&=
\sigma^{\mathrm{T}_{\Sigma^{\boldsymbol{\mathcal{A}}}}(X)}
\left(\left(
\mathrm{CH}^{(1)}_{s_{j}}\left(
\mathfrak{P}'_{j}
\circ_{s_{j}}^{0\mathbf{Pth}_{\boldsymbol{\mathcal{A}}}}
\mathfrak{P}_{j}
\right)
\right)_{j\in\bb{\mathbf{s}}}
\right)
\\&=
\sigma^{\mathrm{T}_{\Sigma^{\boldsymbol{\mathcal{A}}}}(X)}
\left(\left(
\mathrm{CH}^{(1)}_{s_{j}}\left(
\mathfrak{Q}'_{j}
\circ_{s_{j}}^{0\mathbf{Pth}_{\boldsymbol{\mathcal{A}}}}
\mathfrak{Q}_{j}
\right)
\right)_{j\in\bb{\mathbf{s}}}
\right)
\\&=
\sigma^{\mathrm{T}_{\Sigma^{\boldsymbol{\mathcal{A}}}}(X)}
\left(\left(
\mathrm{CH}^{(1)}_{s_{j}}\left(\left(
\mathfrak{Q}'
\circ_{s}^{0\mathbf{Pth}_{\boldsymbol{\mathcal{A}}}}
\mathfrak{Q}
\right)_{j}
\right)\right)_{j\in\bb{\mathbf{s}}}
\right)
\\&=
\mathrm{CH}^{(1)}_{s}\left(
\mathfrak{Q}'
\circ_{s}^{0\mathbf{Pth}_{\boldsymbol{\mathcal{A}}}}
\mathfrak{Q}
\right).
\end{align*}

This completes the case (2.3).

This finishes the proof.
\end{proof}

\chapter{
\texorpdfstring
{On the quotient $[\mathrm{Pth}_{\boldsymbol{\mathcal{A}}}]$}
{On the quotient of paths}
}\label{S1F}

This chapter is devoted to the study of the quotient $\mathrm{Pth}_{\boldsymbol{\mathcal{A}}}/{\mathrm{Ker}(\mathrm{CH}^{(1)})}$, which will be denoted by $[\mathrm{Pth}_{\boldsymbol{\mathcal{A}}}]$ for short. 
In the first place, since the equivalence ${\mathrm{Ker}(\mathrm{CH}^{(1)})}$ on $\mathrm{Pth}_{\boldsymbol{\mathcal{A}}}$ is a closed $\Sigma^{\boldsymbol{\mathcal{A}}}$-congruence on the partial $\Sigma^{\boldsymbol{\mathcal{A}}}$-algebra $\mathbf{Pth}_{\boldsymbol{\mathcal{A}}}$, we equip $[\mathrm{Pth}_{\boldsymbol{\mathcal{A}}}]$ with a structure of partial $\Sigma^{\boldsymbol{\mathcal{A}}}$-algebra, we denote by $[\mathbf{Pth}_{\boldsymbol{\mathcal{A}}}]$ the corresponding partial $\Sigma^{\boldsymbol{\mathcal{A}}}$-algebra. We then define $\Sigma$-homomorphisms from an to $[\mathbf{Pth}_{\boldsymbol{\mathcal{A}}}]$. We show that $[\mathrm{Pth}_{\boldsymbol{\mathcal{A}}}]$ is equipped with a structure of categorial $\Sigma$-algebra, we denote by $[\mathsf{Pth}_{\boldsymbol{\mathcal{A}}}]$ the corresponding categorial $\Sigma$-algebra. We conclude this chapter by showing that $\coprod[\mathrm{Pth}_{\boldsymbol{\mathcal{A}}}]$ is equipped with a structure of Artinian order $\leq_{[\mathbf{Pth}_{\boldsymbol{\mathcal{A}}}]}$, we denote by $(\coprod[\mathrm{Pth}_{\boldsymbol{\mathcal{A}}}],\leq_{[\mathbf{Pth}_{\boldsymbol{\mathcal{A}}}]})$ the corresponding Artinian ordered set. In addition, we prove that several $S$-sorted mappings to $(\coprod[\mathrm{Pth}_{\boldsymbol{\mathcal{A}}}],\leq_{[\mathbf{Pth}_{\boldsymbol{\mathcal{A}}}]})$ are order preseving and some of them are, in addition, order-reflecting or embeddings.


\begin{restatable}{convention}{CCHClass}
\label{CCHClass} 
\index{path!first-order!$[\mathfrak{P}]_{s}$}
\index{path!first-order!$[\mathrm{Pth}_{\boldsymbol{\mathcal{A}}}]$}
To simplify the notation, for a sort $s\in S$ and a path $\mathfrak{P}\in\mathrm{Pth}_{\boldsymbol{\mathcal{A}},s}$, we will let $[\mathfrak{P}]_{s}$ stand for $[\mathfrak{P}]_{\mathrm{Ker}(\mathrm{CH}^{(1)})_{s}}$, the 
$\mathrm{Ker}(\mathrm{CH}^{(1)})_{s}$-equivalence class of $\mathfrak{P}$, and we will call it the \emph{path class of} $\mathfrak{P}$. Moreover, the $S$-sorted quotient $\mathrm{Pth}_{\boldsymbol{\mathcal{A}}}/{\mathrm{Ker}(\mathrm{CH})}$ will simply be denoted by $[\mathrm{Pth}_{\boldsymbol{\mathcal{A}}}]$.
\end{restatable}

\begin{restatable}{definition}{DCHQuot}
\index{projection!first-order!$\mathrm{pr}^{\mathrm{Ker}(\mathrm{CH}^{(1)})}$}
\index{Curry-Howard!first-order!$\mathrm{CH}^{(1)\mathrm{m}}$}
\label{DCHQuot} The (epi,mono)-factorization of the $S$-sorted mapping $\mathrm{CH}^{(1)}$ from 
$\mathrm{Pth}_{\boldsymbol{\mathcal{A}}}$ to $\mathrm{T}_{\Sigma^{\boldsymbol{\mathcal{A}}}}(X)$ is given by
\begin{enumerate}
\item the projection $\mathrm{pr}^{\mathrm{Ker}(\mathrm{CH}^{(1)})}$ from $\mathrm{Pth}_{\boldsymbol{\mathcal{A}}}$ to $[\mathrm{Pth}_{\boldsymbol{\mathcal{A}}}]$ that, for every sort $s\in S$, maps a  path $\mathfrak{P}$ in $\mathrm{Pth}_{\boldsymbol{\mathcal{A}},s}$ to $[\mathfrak{P}]_{s}$, its equivalence class under the kernel of the  Curry-Howard mapping; and
\item the embedding $\mathrm{CH}^{(1)\mathrm{m}}$ from $[\mathrm{Pth}_{\boldsymbol{\mathcal{A}}}]$ to $\mathrm{T}_{\Sigma^{\boldsymbol{\mathcal{A}}}}(X)$ that, for every sort $s\in S$, assigns to an equivalence class $[\mathfrak{P}]_{s}$ in $[\mathrm{Pth}_{\boldsymbol{\mathcal{A}}}]_{s}$ the term $\mathrm{CH}^{(1)}_{s}(\mathfrak{P})$, i.e., the value of the  Curry-Howard mapping at any equivalence class representative. We will refer to the mapping $\mathrm{CH}^{(1)\mathrm{m}}$ as the \emph{monomorphic Curry-Howard mapping}.
\end{enumerate}

The reader is advised to consult the diagram appearing in Figure~\ref{FCHQuot}.
\end{restatable}

\begin{figure}
\begin{center}
\begin{tikzpicture}
[ACliment/.style={-{To [angle'=45, length=5.75pt, width=4pt, round]}},scale=.8]
\node[] (p) at (0,0) [] {$\mathrm{Pth}_{\boldsymbol{\mathcal{A}}}$};
\node[] (pq) at (6,0) [] {$[\mathrm{Pth}_{\boldsymbol{\mathcal{A}}}]
$};
\node[] (t) at (0,-3) [] {$\mathrm{T}_{\Sigma^{\boldsymbol{\mathcal{A}}}}(X)$};

\draw[ACliment]  (p) 	to node [above]	
{$\mathrm{pr}^{\mathrm{Ker}(\mathrm{CH}^{(1)})}$} (pq);
\draw[ACliment]  (p) 	to node [left]	
{$\mathrm{CH}^{(1)}$} (t);
\draw[ACliment, bend left=10]  (pq) 	to node [below right]	
{$\mathrm{CH}^{(1)\mathrm{m}}$} (t);

\end{tikzpicture}
\end{center}
\caption{Many-sorted quotient mappings at layer 1.}
\label{FCHQuot}
\end{figure}

We next define several mappings from and to $[\mathrm{Pth}_{\boldsymbol{\mathcal{A}}}]$.

\begin{restatable}{definition}{DCHEch}
\label{DCHEch} We will denote by
\begin{enumerate}
\item $\mathrm{ip}^{([1],X)}$ the $S$-sorted mapping from $X$ to $[\mathrm{Pth}_{\boldsymbol{\mathcal{A}}}]$ given by the composition $\mathrm{ip}^{([1],X)}=\mathrm{pr}^{\mathrm{Ker}(\mathrm{CH}^{(1)})}\circ \mathrm{ip}^{(1,X)}$, i.e., for every sort $s\in S$, $\mathrm{ip}^{([1],X)}_{s}$ sends a variable $x\in X_{s}$ to the class $[\mathrm{ip}^{(1,X)}_{s}(x)]_{s}$ in $[\mathrm{Pth}_{\boldsymbol{\mathcal{A}}}]_{s}$.
\index{identity!first-order!$\mathrm{ip}^{([1],X)}$}
\item $\mathrm{ech}^{([1],\mathcal{A})}$ the $S$-sorted mapping from $\mathcal{A}$ to $[\mathrm{Pth}_{\boldsymbol{\mathcal{A}}}]$ given by the composition $\mathrm{ech}^{([1],\mathcal{A})}=\mathrm{pr}^{\mathrm{Ker}(\mathrm{CH}^{(1)})}\circ\mathrm{ech}^{(1,\mathcal{A})}$, i.e., for every sort $s\in S$, $\mathrm{ech}^{([1],\mathcal{A})}_{s}$ sends a rewrite rule $\mathfrak{p}\in\mathcal{A}_{s}$ to the class $[\mathrm{ech}^{(1,\mathcal{A})}_{s}(\mathfrak{p})]_{s}$ in $[\mathrm{Pth}_{\boldsymbol{\mathcal{A}}}]_{s}$.
\index{echelon!first-order!$\mathrm{ech}^{([1],\mathcal{A})}$}
\end{enumerate}

The above $S$-sorted mappings are depicted in the diagram of Figure~\ref{FCHEch}.
\end{restatable}

\begin{figure}
\begin{tikzpicture}
[ACliment/.style={-{To [angle'=45, length=5.75pt, width=4pt, round]}},scale=.8]
\node[] (x) at (0,0) [] {$X$};
\node[] (a) at (0,-1.5) [] {$\mathcal{A}$};
\node[] (T) at (6,-1.5) [] {$
[\mathrm{Pth}_{\boldsymbol{\mathcal{A}}}]$};
\draw[ACliment, bend left=10]  (x) to node [above right] {$\mathrm{ip}^{([1],X)}$} (T);
\draw[ACliment]  (a) to node [below] {$\mathrm{ech}^{([1],\mathcal{A})}$} (T);
\end{tikzpicture}
\caption{Quotient path embeddings relative to $X$ and $\mathcal{A}$.}\label{FCHEch}
\end{figure}


\begin{restatable}{definition}{DCHUZ}
\label{DCHUZ} We will denote by

\begin{enumerate}
\item $\mathrm{sc}^{(0,[1])}$ the $S$-sorted mapping from $[\mathrm{Pth}_{\boldsymbol{\mathcal{A}}}]$ to $\mathrm{T}_{\Sigma}(X)$ that, for every $s\in S$, assigns to an equivalence class $[\mathfrak{P}]_{s}$ in $[\mathrm{Pth}_{\boldsymbol{\mathcal{A}}}]_{s}$ with $\mathfrak{P}\in\mathrm{Pth}_{\boldsymbol{\mathcal{A}},s}$ the  term $\mathrm{sc}^{(0,1)}_{s}(\mathfrak{P})$, i.e., the value of the mapping $\mathrm{sc}^{(0,1)}_{s}$ at any equivalence class representative. That $\mathrm{sc}^{(0,[1])}$ is well-defined follows from Lemma~\ref{LCH};
\index{source!first-order!$\mathrm{sc}^{(0,[1])}$}
\item $\mathrm{tg}^{(0,[1])}$ the $S$-sorted mapping from $[\mathrm{Pth}_{\boldsymbol{\mathcal{A}}}]$ to $\mathrm{T}_{\Sigma}(X)$ that, for every $s\in S$, assigns to an equivalence class $[\mathfrak{P}]_{s}$ in $[\mathrm{Pth}_{\boldsymbol{\mathcal{A}}}]_{s}$ with $\mathfrak{P}\in\mathrm{Pth}_{\boldsymbol{\mathcal{A}},s}$ the  term $\mathrm{tg}^{(0,1)}_{s}(\mathfrak{P})$, i.e., the value of the mapping $\mathrm{tg}^{(0,1)}_{s}$ at any equivalence class representative. That $\mathrm{tg}^{(0,[1])}$ is well-defined follows from Lemma~\ref{LCH};
\index{target!first-order!$\mathrm{tg}^{(0,[1])}$}
\item  $\mathrm{ip}^{([1],0)\sharp}$  the $S$-sorted mapping from $\mathrm{T}_{\Sigma}(X)$ to $[\mathrm{Pth}_{\boldsymbol{\mathcal{A}}}]$ given by the composition
$\mathrm{ip}^{([1],0)\sharp}=\mathrm{pr}^{\mathrm{Ker}(\mathrm{CH}^{(1)})}\circ\mathrm{ip}^{(1,0)\sharp}$. That is, for every sort $s\in S$, $\mathrm{ip}^{([1],0)\sharp}$ assigns to a term $P$ in $\mathrm{T}_{\Sigma}(X)_{s}$ the class with respect to the  Curry-Howard mapping of the $(1,0)$-identity  path $\mathrm{ip}^{(1,0)\sharp}_{s}(P)$, i.e., 
$
\mathrm{ip}^{([1],0)\sharp}_{s}(
P
)=[\mathrm{ip}^{(1,0)\sharp}_{s}(P)]_{s}.
$
\index{identity!first-order!$\mathrm{ip}^{([1],0)\sharp}$}
\end{enumerate}

The above $S$-sorted mappings are depicted in the diagram of Figure~\ref{FCHUZ}.
\end{restatable}

\begin{figure}
\begin{center}
\begin{tikzpicture}
[ACliment/.style={-{To [angle'=45, length=5.75pt, width=4pt, round]}},scale=1]
\node[] (x) at (0,0) [] {$X$};
\node[] (t) at (6,0) [] {$\mathrm{T}_{\Sigma}(X)$};
\node[] (pq) at (6,-3) []  {$[\mathrm{Pth}_{\boldsymbol{\mathcal{A}}}]
$};

\draw[ACliment]  (x) 	to node [above right]	
{$\eta^{(0,X)}$} (t);
\draw[ACliment, bend right=10]  (x) 	to node [below left]	
{$\mathrm{ip}^{([1],X)}$} (pq);

\node[] (B1) at (6,-1.5)  [] {};
\draw[ACliment]  ($(B1)+(0,1.2)$) to node [above, fill=white] {
$\mathrm{ip}^{([1],0)\sharp}$
} ($(B1)+(0,-1.2)$);
\draw[ACliment, bend right=40]  ($(B1)+(.3,-1.2)$) to node [ below, fill=white] {
$\mathrm{tg}^{(0,[1])}$
} ($(B1)+(.3,1.2)$);
\draw[ACliment, bend left=40]  ($(B1)+(-.3,-1.2)$) to node [below, fill=white] {
$ \mathrm{sc}^{(0,[1])}$
} ($(B1)+(-.3,1.2)$);

\end{tikzpicture}
\end{center}
\caption{Many-sorted quotient mappings relative to $X$ at layers 0 \& 1.}
\label{FCHUZ}
\end{figure}

The following are some of the relationships between the $S$-sorted mappings $\eta^{(0,X)}$, $\mathrm{ip}^{([1],X)}$, $\mathrm{sc}^{(0,[1])}$, $\mathrm{tg}^{(0,[1])}$, $\mathrm{ip}^{([1],0)\sharp}$ and 
$\mathrm{id}^{\mathrm{T}_{\Sigma}(X)}$. 

\begin{proposition}\label{PCHBasicEq} The following equalities holds
\begin{itemize}
\item[(i)] $\mathrm{sc}^{(0,[1])}\circ\mathrm{ip}^{([1],0)\sharp}=\mathrm{id}^{\mathrm{T}_{\Sigma}(X)}$;
\item[(ii)] $\mathrm{tg}^{(0,[1])}\circ\mathrm{ip}^{([1],0)\sharp}=\mathrm{id}^{\mathrm{T}_{\Sigma}(X)}$.
\end{itemize}
\end{proposition}

\begin{proposition} The following equalities holds
\begin{itemize}
\item[(i)] $\mathrm{sc}^{(0,[1])}\circ\mathrm{ip}^{([1],X)}=\eta^{(0,X)};$
\item[(ii)] $\mathrm{tg}^{(0,[1])}\circ\mathrm{ip}^{([1],X)}=\eta^{(0,X)};$
\item[(iii)] $\mathrm{ip}^{([1],0)\sharp}\circ\eta^{(0,X)}=\mathrm{ip}^{([1],X)}$.
\end{itemize}

The reader is advised to consult the diagram appearing in Figure~\ref{FCHUZ}.
\end{proposition}

\section{
\texorpdfstring
{A structure of partial  $\Sigma^{\boldsymbol{\mathcal{A}}}$-algebra on $[\mathrm{Pth}_{\boldsymbol{\mathcal{A}}}]$}
{An algebra on the quotient of paths}
}

As an immediate consequence of the main result of the above section, we next state that the $S$-sorted set $[\mathrm{Pth}_{\boldsymbol{\mathcal{A}}}]$ is equipped with a structure of partial $\Sigma^{\boldsymbol{\mathcal{A}}}$-algebra. 

\begin{restatable}{proposition}{PCHCatAlg}
\label{PCHCatAlg}
\index{path!first-order!$[\mathbf{Pth}_{\boldsymbol{\mathcal{A}}}]$}
The $S$-sorted set $[\mathrm{Pth}_{\boldsymbol{\mathcal{A}}}]$ is equipped, in a natural way, with a structure of  partial $\Sigma^{\boldsymbol{\mathcal{A}}}$-algebra, we denote by $[\mathbf{Pth}_{\boldsymbol{\mathcal{A}}}]$ the corresponding $\Sigma^{\boldsymbol{\mathcal{A}}}$-algebra (which is a quotient of $\mathbf{Pth}_{\boldsymbol{\mathcal{A}}}$, the partial $\Sigma^{\boldsymbol{\mathcal{A}}}$-algebra constructed in Proposition~\ref{PPthCatAlg}). For later reference, we point out that the operations of $[\mathbf{Pth}_{\boldsymbol{\mathcal{A}}}]$ are denoted by $\mathrm{sc}^{0[\mathbf{Pth}_{\boldsymbol{\mathcal{A}}}]}$, 
$\mathrm{tg}^{0[\mathbf{Pth}_{\boldsymbol{\mathcal{A}}}]}$, $\circ^{0[\mathbf{Pth}_{\boldsymbol{\mathcal{A}}}]}$ and 
$\sigma^{[\mathbf{Pth}_{\boldsymbol{\mathcal{A}}}]}$. Moreover, the $S$-sorted mapping 
$$
\mathrm{pr}^{\mathrm{Ker}(\mathrm{CH}^{(1)})}
\colon
\mathbf{Pth}_{\boldsymbol{\mathcal{A}}}
\mor
[\mathbf{Pth}_{\boldsymbol{\mathcal{A}}}]
$$
is a closed and surjective $\Sigma^{\boldsymbol{\mathcal{A}}}$-homomorphism from $\mathbf{Pth}_{\boldsymbol{\mathcal{A}}}$ to $[\mathbf{Pth}_{\boldsymbol{\mathcal{A}}}]$.
\end{restatable}

\begin{remark}\label{RCHMono} The monomorphic Curry-Howard mapping, that is,
$$
\mathrm{CH}^{(1)\mathrm{m}}
\colon
[\mathrm{Pth}_{\boldsymbol{\mathcal{A}}}]
\mor
\mathrm{T}_{\Sigma^{\boldsymbol{\mathcal{A}}}}(X)
$$
is not a $\Sigma^{\boldsymbol{\mathcal{A}}}$-homomorphism from 
$[\mathbf{Pth}_{\boldsymbol{\mathcal{A}}}]$ to 
$\mathbf{T}_{\Sigma^{\boldsymbol{\mathcal{A}}}}(X)$. The reader is advised to refer to the counterexample in Proposition~\ref{PCHNotHomCat}. Later on, by considering a suitable congruence on $\mathbf{T}_{\Sigma^{\boldsymbol{\mathcal{A}}}}(X)$ and composing $\mathrm{CH}^{(1)\mathrm{m}}$ with the canonical projection to the quotient we will obtain a $\Sigma^{\boldsymbol{\mathcal{A}}}$-homomorphism.
\end{remark}

We next consider the $\Sigma$-reduct of the partial $\Sigma^{\boldsymbol{\mathcal{A}}}$-algebra 
$[\mathbf{Pth}_{\boldsymbol{\mathcal{A}}}]$ and investigate its connections with the mappings defined in the previous section.

\begin{definition} For the partial $\Sigma^{\boldsymbol{\mathcal{A}}}$-algebra $[\mathbf{Pth}_{\boldsymbol{\mathcal{A}}}]$, we denote by $[\mathbf{Pth}^{(0,1)}_{\boldsymbol{\mathcal{A}}}]$ the $\Sigma$-algebra 
$\mathbf{in}^{\Sigma,(0,1)}\left([\mathbf{Pth}_{\boldsymbol{\mathcal{A}}}]\right)$. We will call  
$[\mathbf{Pth}^{(0,1)}_{\boldsymbol{\mathcal{A}}}]$ the $\Sigma$-reduct of the partial 
$\Sigma^{\boldsymbol{\mathcal{A}}}$-algebra $[\mathbf{Pth}_{\boldsymbol{\mathcal{A}}}]$.
\end{definition}

\begin{proposition} The $S$-sorted mapping $\mathrm{pr}^{\mathrm{Ker}(\mathrm{CH}^{(1)})}$
is a closed and surjective $\Sigma$-homomorphism from $\mathbf{Pth}^{(0,1)}_{\boldsymbol{\mathcal{A}}}$ to $[\mathbf{Pth}^{(0,1)}_{\boldsymbol{\mathcal{A}}}]$.
\end{proposition}

\begin{restatable}{proposition}{PCHDZ}
\label{PCHDZ} The $S$-sorted  mappings $\mathrm{sc}^{(0,[1])}$ and $\mathrm{tg}^{(0,[1])}$ are $\Sigma$-homomorphisms from $[\mathbf{Pth}^{(0,1)}_{\boldsymbol{\mathcal{A}}}]$ to $\mathbf{T}_{\Sigma}(X)$.
\end{restatable}

\begin{restatable}{proposition}{PCHDZIp}
\label{PCHDZIp} The $S$-sorted mapping $\mathrm{ip}^{([1],0)\sharp}$  is a $\Sigma$-homomorphism from  
$\mathbf{T}_{\Sigma}(X)$ to $[\mathbf{Pth}^{(0,1)}_{\boldsymbol{\mathcal{A}}}]$.
\end{restatable}

\section{
\texorpdfstring
{A structure of $S$-sorted categorial $\Sigma$-algebra on $[\mathrm{Pth}_{\boldsymbol{\mathcal{A}}}]$}
{A categorial algebra on the quotient of paths}
}

In what follows we will prove that the $S$-sorted set $[\mathrm{Pth}_{\boldsymbol{\mathcal{A}}}]$ is equipped with a structure of $S$-sorted categorial $\Sigma$-algebra. To do so we begin by showing that the partial 
$\Sigma^{\boldsymbol{\mathcal{A}}}$-algebra $[\mathbf{Pth}_{\boldsymbol{\mathcal{A}}}]$ satisfies the defining equations of the concept of $S$-sorted category.

\begin{proposition}\label{PCHVarA2} Let $s$ be a sort in $S$ and $[\mathfrak{P}]_{s}$ a path class in $[\mathrm{Pth}_{\boldsymbol{\mathcal{A}}}]_{s}$, then the following equalities holds
\begin{align*}
\mathrm{sc}^{0[\mathbf{Pth}_{\boldsymbol{\mathcal{A}}}]}_{s}\left(
\mathrm{sc}^{0[\mathbf{Pth}_{\boldsymbol{\mathcal{A}}}]}_{s}\left(
[\mathfrak{P}]_{s}
\right)\right)
&=
\mathrm{sc}^{0[\mathbf{Pth}_{\boldsymbol{\mathcal{A}}}]}_{s}\left(
[\mathfrak{P}]_{s}
\right);
\\
\mathrm{sc}^{0[\mathbf{Pth}_{\boldsymbol{\mathcal{A}}}]}_{s}\left(
\mathrm{tg}^{0[\mathbf{Pth}_{\boldsymbol{\mathcal{A}}}]}_{s}\left(
[\mathfrak{P}]_{s}
\right)\right)
&=
\mathrm{tg}^{0[\mathbf{Pth}_{\boldsymbol{\mathcal{A}}}]}_{s}\left(
[\mathfrak{P}]_{s}
\right);
\\
\mathrm{tg}^{0[\mathbf{Pth}_{\boldsymbol{\mathcal{A}}}]}_{s}\left(
\mathrm{sc}^{0[\mathbf{Pth}_{\boldsymbol{\mathcal{A}}}]}_{s}\left(
[\mathfrak{P}]_{s}
\right)\right)
&=
\mathrm{sc}^{0[\mathbf{Pth}_{\boldsymbol{\mathcal{A}}}]}_{s}\left(
[\mathfrak{P}]_{s}
\right);
\\
\mathrm{tg}^{0[\mathbf{Pth}_{\boldsymbol{\mathcal{A}}}]}_{s}\left(
\mathrm{tg}^{0[\mathbf{Pth}_{\boldsymbol{\mathcal{A}}}]}_{s}\left(
[\mathfrak{P}]_{s}
\right)\right)
&=
\mathrm{tg}^{0[\mathbf{Pth}_{\boldsymbol{\mathcal{A}}}]}_{s}\left(
[\mathfrak{P}]_{s}
\right).
\end{align*}
\end{proposition}
\begin{proof}
Regarding the first equation, the following chain of equalities holds
\allowdisplaybreaks
\begin{align*}
\mathrm{sc}^{0[\mathbf{Pth}_{\boldsymbol{\mathcal{A}}}]}_{s}\left(
\mathrm{sc}^{0[\mathbf{Pth}_{\boldsymbol{\mathcal{A}}}]}_{s}\left(
[\mathfrak{P}]_{s}
\right)\right)
&=
\left[
\mathrm{sc}^{0\mathbf{Pth}_{\boldsymbol{\mathcal{A}}}}_{s}\left(
\mathrm{sc}^{0\mathbf{Pth}_{\boldsymbol{\mathcal{A}}}}_{s}\left(
\mathfrak{P}
\right)\right)\right]_{s}
\tag{1}
\\&=
\left[
\mathrm{ip}^{(1,0)\sharp}_{s}\left(
\mathrm{sc}^{(0,1)}_{s}\left(
\mathrm{ip}^{(1,0)\sharp}_{s}\left(
\mathrm{sc}^{(0,1)}_{s}\left(
\mathfrak{P}
\right)\right)\right)\right)
\right]_{s}
\tag{2}
\\&=
\left[
\mathrm{ip}^{(1,0)\sharp}_{s}\left(
\mathrm{sc}^{(0,1)}_{s}\left(
\mathfrak{P}
\right)\right)
\right]_{s}
\tag{3}
\\&=
\left[
\mathrm{sc}^{0\mathbf{Pth}_{\boldsymbol{\mathcal{A}}}}_{s}\left(
\mathfrak{P}
\right)
\right]_{s}
\tag{4}
\\&=
\mathrm{sc}^{0[\mathbf{Pth}_{\boldsymbol{\mathcal{A}}}]}_{s}\left(
[\mathfrak{P}]_{s}
\right).
\tag{5}
\end{align*}

The first equality unravels the interpretation of $0$-source operation in the many-sorted partial $\Sigma^{\boldsymbol{\mathcal{A}}}$-algebra $[\mathbf{Pth}_{\boldsymbol{\mathcal{A}}}]$ according to Proposition~\ref{PCHCatAlg}; the second equality unravels the interpretation of the $0$-source operation in the many-sorted partial $\Sigma^{\boldsymbol{\mathcal{A}}}$-algebra $\mathbf{Pth}_{\boldsymbol{\mathcal{A}}}$ according to Proposition~\ref{PPthCatAlg}; the third equality follows from Proposition~\ref{PBasicEq}; the fourth equality recovers the interpretation of the $0$-source operation in the many-sorted partial $\Sigma^{\boldsymbol{\mathcal{A}}}$-algebra $\mathbf{Pth}_{\boldsymbol{\mathcal{A}}}$ according to Proposition~\ref{PPthCatAlg}; finally, the last equality recovers the interpretation of $0$-source operation in the many-sorted partial $\Sigma^{\boldsymbol{\mathcal{A}}}$-algebra $[\mathbf{Pth}_{\boldsymbol{\mathcal{A}}}]$ according to Proposition~\ref{PCHCatAlg}.

The remaining equalities are handled similarly.

This completes the proof.
\end{proof}

\begin{proposition}\label{PCHVarA3} Let $s$ be a sort in $S$ and $[\mathfrak{P}]_{s}$, $[\mathfrak{Q}]_{s}$ path classes in $[\mathrm{Pth}_{\boldsymbol{\mathcal{A}}}]_{s}$, then the following statements are equivalent
\begin{enumerate}
\item[(i)] $[\mathfrak{Q}]_{s}
\circ^{0[\mathbf{Pth}_{\boldsymbol{\mathcal{A}}}]}_{s}
[\mathfrak{P}]_{s}$ is defined;
\item[(ii)] $\mathrm{sc}^{0[\mathbf{Pth}_{\boldsymbol{\mathcal{A}}}]}_{s}\left(
[\mathfrak{Q}]_{s}
\right)
=
\mathrm{tg}^{0[\mathbf{Pth}_{\boldsymbol{\mathcal{A}}}]}_{s}\left(
[\mathfrak{P}]_{s}
\right).$
\end{enumerate}
\end{proposition}
\begin{proof}
The following chain of equivalences holds
\begin{flushleft}
$[\mathfrak{Q}]_{s}
\circ^{0[\mathbf{Pth}_{\boldsymbol{\mathcal{A}}}]}_{s}
[\mathfrak{P}]_{s}\mbox{ is defined }$
\allowdisplaybreaks
\begin{align*}
\qquad
&\Leftrightarrow
\mathfrak{Q}\circ^{0\mathbf{Pth}_{\boldsymbol{\mathcal{A}}}}_{s}\mathfrak{P}
\mbox{ is defined }
\tag{1}
\\&\Leftrightarrow
\mathrm{sc}^{(0,1)}_{s}\left(
\mathfrak{Q}
\right)
=
\mathrm{tg}^{(0,1)}_{s}\left(
\mathfrak{P}
\right)
\tag{2}
\\&\Leftrightarrow
\mathrm{ip}^{(1,0)\sharp}_{s}\left(
\mathrm{sc}^{(0,1)}_{s}\left(
\mathfrak{Q}
\right)\right)
=
\mathrm{ip}^{(1,0)\sharp}_{s}\left(
\mathrm{tg}^{(0,1)}_{s}\left(
\mathfrak{P}
\right)\right)
\tag{3}
\\&\Leftrightarrow
\left[
\mathrm{ip}^{(1,0)\sharp}_{s}\left(
\mathrm{sc}^{(0,1)}_{s}\left(
\mathfrak{Q}
\right)\right)\right]_{s}
=
\left[
\mathrm{ip}^{(1,0)\sharp}_{s}\left(
\mathrm{tg}^{(0,1)}_{s}\left(
\mathfrak{P}
\right)\right)
\right]_{s}
\tag{4}
\\&\Leftrightarrow
\left[
\mathrm{sc}^{0\mathbf{Pth}_{\boldsymbol{\mathcal{A}}}}_{s}\left(
\mathfrak{Q}
\right)\right]_{s}
=
\left[
\mathrm{tg}^{0\mathbf{Pth}_{\boldsymbol{\mathcal{A}}}}_{s}\left(
\mathfrak{P}
\right)\right]_{s}
\tag{5}
\\&\Leftrightarrow
\mathrm{sc}^{0[\mathbf{Pth}_{\boldsymbol{\mathcal{A}}}]}_{s}\left(
[\mathfrak{Q}]_{s}
\right)
=
\mathrm{tg}^{0[\mathbf{Pth}_{\boldsymbol{\mathcal{A}}}]}_{s}\left(
[\mathfrak{P}]_{s}
\right).
\tag{6}
\end{align*}
\end{flushleft}

In the just stated chain of equivalences, the first equivalence follows from Proposition~\ref{PCHCatAlg}; the second equivalence follows from Definition~\ref{DPthComp}; the third equivalence follows from Proposition~\ref{PBasicEq}; the fourth equivalence follows from Corollary~\ref{CCHUZId}; the fifth equivalence follows from the description of the $0$-source and $0$-target operation symbols in the many-sorted partial $\Sigma^{\boldsymbol{\mathcal{A}}}$-algebra according to Proposition~\ref{PPthCatAlg}; finally, the last equivalence follows from the description of the $0$-source and $0$-target operation symbols in the many-sorted partial $\Sigma^{\boldsymbol{\mathcal{A}}}$-algebra $[\mathbf{Pth}_{\boldsymbol{\mathcal{A}}}]$, according to Proposition~\ref{PCHCatAlg}.

This completes the proof.
\end{proof}

\begin{proposition}\label{PCHVarA4} Let $s$ be a sort in $S$ and $[\mathfrak{P}]_{s}$, $[\mathfrak{Q}]_{s}$ path classes in $[\mathrm{Pth}_{\boldsymbol{\mathcal{A}}}]_{s}$. If 
$\mathrm{sc}^{0[\mathbf{Pth}_{\boldsymbol{\mathcal{A}}}]}_{s}\left(
[\mathfrak{Q}]_{s}
\right)
=
\mathrm{tg}^{0[\mathbf{Pth}_{\boldsymbol{\mathcal{A}}}]}_{s}\left(
[\mathfrak{P}]_{s}
\right)$, then the following equalities hold
\allowdisplaybreaks
\begin{align*}
\mathrm{sc}^{0[\mathbf{Pth}_{\boldsymbol{\mathcal{A}}}]}_{s}\left(
[\mathfrak{Q}]_{s}
\circ^{0[\mathbf{Pth}_{\boldsymbol{\mathcal{A}}}]}_{s}
[\mathfrak{P}]_{s}
\right)
&=
\mathrm{sc}^{0[\mathbf{Pth}_{\boldsymbol{\mathcal{A}}}]}_{s}\left(
[\mathfrak{P}]_{s}
\right);
\\
\mathrm{tg}^{0[\mathbf{Pth}_{\boldsymbol{\mathcal{A}}}]}_{s}\left(
[\mathfrak{Q}]_{s}
\circ^{0[\mathbf{Pth}_{\boldsymbol{\mathcal{A}}}]}_{s}
[\mathfrak{P}]_{s}
\right)
&=
\mathrm{tg}^{0[\mathbf{Pth}_{\boldsymbol{\mathcal{A}}}]}_{s}\left(
[\mathfrak{Q}]_{s}
\right).
\end{align*}
\end{proposition}
\begin{proof}
The following chain of equalities holds
\allowdisplaybreaks
\begin{align*}
\mathrm{sc}^{0[\mathbf{Pth}_{\boldsymbol{\mathcal{A}}}]}_{s}\left(
[\mathfrak{Q}]_{s}
\circ^{0[\mathbf{Pth}_{\boldsymbol{\mathcal{A}}}]}_{s}
[\mathfrak{P}]_{s}
\right)&=
\mathrm{sc}^{0[\mathbf{Pth}_{\boldsymbol{\mathcal{A}}}]}_{s}\left(
\left[\mathfrak{Q}
\circ^{0\mathbf{Pth}_{\boldsymbol{\mathcal{A}}}}_{s}
\mathfrak{P}
\right]_{s}
\right)
\tag{1}
\\&=
\left[
\mathrm{sc}^{0\mathbf{Pth}_{\boldsymbol{\mathcal{A}}}}_{s}\left(
\mathfrak{Q}
\circ^{0\mathbf{Pth}_{\boldsymbol{\mathcal{A}}}}_{s}
\mathfrak{P}
\right)
\right]_{s}
\tag{2}
\\&=
\left[
\mathrm{ip}^{(1,0)\sharp}_{s}\left(
\mathrm{sc}^{(0,1)}_{s}\left(
\mathfrak{Q}
\circ^{0\mathbf{Pth}_{\boldsymbol{\mathcal{A}}}}_{s}
\mathfrak{P}
\right)
\right)
\right]_{s}
\tag{3}
\\&=
\left[
\mathrm{ip}^{(1,0)\sharp}_{s}\left(
\mathrm{sc}^{(0,1)}_{s}\left(
\mathfrak{P}
\right)
\right)
\right]_{s}
\tag{4}
\\&=
\left[
\mathrm{sc}^{0\mathbf{Pth}_{\boldsymbol{\mathcal{A}}}}_{s}\left(
\mathfrak{P}
\right)
\right]_{s}
\tag{5}
\\&=
\mathrm{sc}^{0[\mathbf{Pth}_{\boldsymbol{\mathcal{A}}}]}_{s}\left(
\left[
\mathfrak{P}
\right]_{s}
\right).
\tag{6}
\end{align*}

The first equality applies the $0$-composition operation in the many-sorted partial $\Sigma^{\boldsymbol{\mathcal{A}}}$-algebra $[\mathbf{Pth}_{\boldsymbol{\mathcal{A}}}]$ according to Proposition~\ref{PCHCatAlg}. In this regard, we recall from Proposition~\ref{PCHVarA3}, that this $0$-composition is well-defined by hypothesis;  the second equality applies the $0$-source operation in the many-sorted partial $\Sigma^{\boldsymbol{\mathcal{A}}}$-algebra $[\mathbf{Pth}_{\boldsymbol{\mathcal{A}}}]$ according to Proposition~\ref{PCHCatAlg}; the third equality unravels the interpretation of the $0$-source operation in the many-sorted partial $\Sigma^{\boldsymbol{\mathcal{A}}}$-algebra $\mathbf{Pth}_{\boldsymbol{\mathcal{A}}}$ according to Proposition~\ref{PPthCatAlg}; the fourth equality follows from Proposition~\ref{PPthComp}; the fifth equality recovers the interpretation of the $0$-source operation in the many-sorted partial $\Sigma^{\boldsymbol{\mathcal{A}}}$-algebra $\mathbf{Pth}_{\boldsymbol{\mathcal{A}}}$ according to Proposition~\ref{PPthCatAlg}; finally, the last equality recovers the $0$-source operation in the many-sorted partial $\Sigma^{\boldsymbol{\mathcal{A}}}$-algebra $[\mathbf{Pth}_{\boldsymbol{\mathcal{A}}}]$ according to Proposition~\ref{PCHCatAlg}.

The remaining equality is handled in a similar way.

This completes the proof.
\end{proof}

\begin{proposition}\label{PCHVarA5} Let $s$ be a sort in $S$ and $[\mathfrak{P}]_{s}$ a path class in $[\mathrm{Pth}_{\boldsymbol{\mathcal{A}}}]_{s}$. Then the following equalities hold
\allowdisplaybreaks
\begin{align*}
[\mathfrak{P}]_{s}
\circ^{0[\mathbf{Pth}_{\boldsymbol{\mathcal{A}}}]}_{s}
\left(
\mathrm{sc}^{0[\mathbf{Pth}_{\boldsymbol{\mathcal{A}}}]}_{s}\left(
[\mathfrak{P}]_{s}
\right)
\right)
&=
[\mathfrak{P}]_{s};
\\
\left(\mathrm{tg}^{0[\mathbf{Pth}_{\boldsymbol{\mathcal{A}}}]}_{s}\left(
[\mathfrak{P}]_{s}
\right)
\right)
\circ^{0[\mathbf{Pth}_{\boldsymbol{\mathcal{A}}}]}_{s}
[\mathfrak{P}]_{s}
&=
[\mathfrak{P}]_{s}.
\end{align*}
\end{proposition}
\begin{proof}
The following chain of equalities holds
\allowdisplaybreaks
\begin{align*}
[\mathfrak{P}]_{s}
\circ^{0[\mathbf{Pth}_{\boldsymbol{\mathcal{A}}}]}_{s}
\left(
\mathrm{sc}^{0[\mathbf{Pth}_{\boldsymbol{\mathcal{A}}}]}_{s}\left(
[\mathfrak{P}]_{s}
\right)
\right)&=
[\mathfrak{P}]_{s}
\circ^{0[\mathbf{Pth}_{\boldsymbol{\mathcal{A}}}]}_{s}
\left(
\left[
\mathrm{sc}^{0\mathbf{Pth}_{\boldsymbol{\mathcal{A}}}}_{s}\left(
\mathfrak{P}
\right)
\right]_{s}
\right)
\tag{1}
\\&=
\left[
\mathfrak{P}
\circ^{0\mathbf{Pth}_{\boldsymbol{\mathcal{A}}}}_{s}
\left(
\mathrm{sc}^{0\mathbf{Pth}_{\boldsymbol{\mathcal{A}}}}_{s}\left(
\mathfrak{P}
\right)
\right)
\right]_{s}
\tag{2}
\\&=
\left[
\mathfrak{P}
\circ^{0\mathbf{Pth}_{\boldsymbol{\mathcal{A}}}}_{s}
\left(
\mathrm{ip}^{(1,0)\sharp}_{s}\left(
\mathrm{sc}^{(0,1)}_{s}\left(
\mathfrak{P}
\right)
\right)
\right)
\right]_{s}
\tag{3}
\\&=
\left[
\mathfrak{P}
\right]_{s}.
\tag{4}
\end{align*}

The first equality unravels the $0$-source operation in the many-sorted partial $\Sigma^{\boldsymbol{\mathcal{A}}}$-algebra $[\mathbf{Pth}_{\boldsymbol{\mathcal{A}}}]$ according to Proposition~\ref{PCHCatAlg}; the second equality unravels the $0$-composition operation in the many-sorted partial $\Sigma^{\boldsymbol{\mathcal{A}}}$-algebra $[\mathbf{Pth}_{\boldsymbol{\mathcal{A}}}]$ according to Proposition~\ref{PCHCatAlg}. In this regard, note that this $0$-composition is defined according to Propositions~\ref{PCHVarA2} and~\ref{PCHVarA3}; the third equality recovers the $0$-source operation in the many-sorted partial $\Sigma^{\boldsymbol{\mathcal{A}}}$-algebra $\mathbf{Pth}_{\boldsymbol{\mathcal{A}}}$ according to Proposition~\ref{PPthCatAlg}; finally, the last equality follows from Proposition~\ref{PPthComp}.

The remaining equality is handled in a similar way.

This completes the proof.
\end{proof}

\begin{proposition}\label{PCHVarA6} Let $s$ be a sort in $S$ and $[\mathfrak{P}]_{s}$, $[\mathfrak{Q}]_{s}$ and $[\mathfrak{R}]_{s}$ path classes in $[\mathrm{Pth}_{\boldsymbol{\mathcal{A}}}]_{s}$ such that 
\allowdisplaybreaks
\begin{align*}
\mathrm{sc}^{0[\mathbf{Pth}_{\boldsymbol{\mathcal{A}}}]}_{s}\left(
[\mathfrak{R}]_{s}
\right)&=
\mathrm{tg}^{0[\mathbf{Pth}_{\boldsymbol{\mathcal{A}}}]}_{s}\left(
[\mathfrak{Q}]_{s}
\right);
\\
\mathrm{sc}^{0[\mathbf{Pth}_{\boldsymbol{\mathcal{A}}}]}_{s}\left(
[\mathfrak{Q}]_{s}
\right)&=
\mathrm{tg}^{0[\mathbf{Pth}_{\boldsymbol{\mathcal{A}}}]}_{s}\left(
[\mathfrak{P}]_{s}
\right). 
\end{align*}
Then the following equality holds
\[
[\mathfrak{R}]_{s}
\circ^{0[\mathbf{Pth}_{\boldsymbol{\mathcal{A}}}]}_{s}\left(
[\mathfrak{Q}]_{s}
\circ^{0[\mathbf{Pth}_{\boldsymbol{\mathcal{A}}}]}_{s}
[\mathfrak{P}]_{s}
\right)
=
\left(
[\mathfrak{R}]_{s}
\circ^{0[\mathbf{Pth}_{\boldsymbol{\mathcal{A}}}]}_{s}
[\mathfrak{Q}]_{s}
\right)
\circ^{0[\mathbf{Pth}_{\boldsymbol{\mathcal{A}}}]}_{s}
[\mathfrak{P}]_{s}
.\]
\end{proposition}
\begin{proof}
The following chain of equalities holds
\allowdisplaybreaks
\begin{align*}
[\mathfrak{R}]_{s}
\circ^{0[\mathbf{Pth}_{\boldsymbol{\mathcal{A}}}]}_{s}\left(
[\mathfrak{Q}]_{s}
\circ^{0[\mathbf{Pth}_{\boldsymbol{\mathcal{A}}}]}_{s}
[\mathfrak{P}]_{s}
\right)&=
\left[
\mathfrak{R}
\circ^{0\mathbf{Pth}_{\boldsymbol{\mathcal{A}}}}_{s}
\left(
\mathfrak{Q}
\circ^{0\mathbf{Pth}_{\boldsymbol{\mathcal{A}}}}_{s}
\mathfrak{P}
\right)
\right]_{s}
\tag{1}
\\&=
\left[
\left(
\mathfrak{R}
\circ^{0\mathbf{Pth}_{\boldsymbol{\mathcal{A}}}}_{s}
\mathfrak{Q}
\right)
\circ^{0\mathbf{Pth}_{\boldsymbol{\mathcal{A}}}}_{s}
\mathfrak{P}s
\right]_{s}
\tag{2}
\\&=
\left(
[\mathfrak{R}]_{s}
\circ^{0[\mathbf{Pth}_{\boldsymbol{\mathcal{A}}}]}_{s}
[\mathfrak{Q}]_{s}
\right)
\circ^{0[\mathbf{Pth}_{\boldsymbol{\mathcal{A}}}]}_{s}
[\mathfrak{P}]_{s}
.
\tag{3}
\end{align*}

The first equality unravels the $0$-composition operation in the many-sorted partial $\Sigma^{\boldsymbol{\mathcal{A}}}$-algebra $[\mathbf{Pth}_{\boldsymbol{\mathcal{A}}}]$ according to Proposition~\ref{PCHCatAlg}. In this regard, note that this $0$-composition is defined according to Proposition~\ref{PCHVarA3}; the second equality follows from Proposition~\ref{PPthComp}; finally, the last equality recovers the $0$-composition operation in the many-sorted partial $\Sigma^{\boldsymbol{\mathcal{A}}}$-algebra $[\mathbf{Pth}_{\boldsymbol{\mathcal{A}}}]$ according to Proposition~\ref{PCHCatAlg}.

This completes the proof.
\end{proof}


\begin{restatable}{proposition}{PCHCat}
\label{PCHCat} 
The ordered pair 
$$
\left(
[\mathrm{Pth}_{\boldsymbol{\mathcal{A}}}],
\left(
\circ^{0[\mathbf{Pth}_{\boldsymbol{\mathcal{A}}}]},
\mathrm{sc}^{0[\mathbf{Pth}_{\boldsymbol{\mathcal{A}}}]},
\mathrm{tg}^{0[\mathbf{Pth}_{\boldsymbol{\mathcal{A}}}]}
\right)
\right),
$$
denoted by $[\mathsf{Pth}_{\boldsymbol{\mathcal{A}}}]$, is an $S$-sorted category.
\end{restatable}

\begin{proof}
Let us recall that, for every $s\in S$, the structure
$([\mathrm{Pth}_{\boldsymbol{\mathcal{A}}}]_{s},\xi_{0,s})$, where
$$
\xi_{0,s}=
\left(
\circ^{0[\mathbf{Pth}_{\boldsymbol{\mathcal{A}}}]}_{s},
\mathrm{sc}^{0[\mathbf{Pth}_{\boldsymbol{\mathcal{A}}}]}_{s},
\mathrm{tg}^{0[\mathbf{Pth}_{\boldsymbol{\mathcal{A}}}]}_{s}
\right)
$$
is a single-sorted category in virtue of Definition~\ref{DCat} and Propositions~\ref{PCHVarA2}, \ref{PCHVarA3}, \ref{PCHVarA4}, \ref{PCHVarA5} and~\ref{PCHVarA6}.
Therefore, by Definition~\ref{DnCatAlg}, the ordered pair
$$
\left(
\left(
[\mathrm{Pth}_{\boldsymbol{\mathcal{A}}}]_{s},\xi_{0,s}
\right)
\right)_{s\in S}
$$
is an $S$-sorted category.
\end{proof}


To complete the proof that $[\mathsf{Pth}_{\boldsymbol{\mathcal{A}}}]$ is a categorial $\Sigma$-algebra, we prove the following propositions.

\begin{proposition}\label{PCHVarA7} Let $(\mathbf{s},s)$ be an element in $S^{\star}\times S$, $\sigma$ an operation symbol in $\Sigma^{\boldsymbol{\mathcal{A}}}_{\mathbf{s},s}$ and $([\mathfrak{P}_{j}]_{s_{j}})_{j\in\bb{\mathbf{s}}}$ a family of path classes in $[\mathrm{Pth}_{\boldsymbol{\mathcal{A}}}]_{\mathbf{s}}$. Then the following equalities holds
\allowdisplaybreaks
\begin{align*}
\sigma^{[\mathbf{Pth}_{\boldsymbol{\mathcal{A}}}]}\left(\left(
\mathrm{sc}^{0[\mathbf{Pth}_{\boldsymbol{\mathcal{A}}}]}_{s_{j}}\left(
[\mathfrak{P}_{j}]_{s_{j}}
\right)
\right)_{j\in\bb{\mathbf{s}}}\right)&=
\mathrm{sc}^{0[\mathbf{Pth}_{\boldsymbol{\mathcal{A}}}]}_{s}\left(
\sigma^{[\mathbf{Pth}_{\boldsymbol{\mathcal{A}}}]}\left(
\left(
[\mathfrak{P}]_{s_{j}}
\right)_{j\in\bb{\mathbf{s}}}
\right)
\right);
\\
\sigma^{[\mathbf{Pth}_{\boldsymbol{\mathcal{A}}}]}\left(\left(
\mathrm{tg}^{0[\mathbf{Pth}_{\boldsymbol{\mathcal{A}}}]}_{s_{j}}\left(
[\mathfrak{P}_{j}]_{s_{j}}
\right)
\right)_{j\in\bb{\mathbf{s}}}\right)&=
\mathrm{tg}^{0[\mathbf{Pth}_{\boldsymbol{\mathcal{A}}}]}_{s}\left(
\sigma^{[\mathbf{Pth}_{\boldsymbol{\mathcal{A}}}]}\left(
\left(
[\mathfrak{P}]_{s_{j}}
\right)_{j\in\bb{\mathbf{s}}}
\right)
\right).
\end{align*}
\end{proposition}
\begin{proof}
We will only present the proof for the first equality. The other one is handled in a similar manner.

Note that the following chain of equalities holds
\allowdisplaybreaks
\begin{align*}
\sigma^{\mathbf{Pth}_{\boldsymbol{\mathcal{A}}}}\left(\left(
\mathrm{sc}^{0\mathbf{Pth}_{\boldsymbol{\mathcal{A}}}}_{s_{j}}\left(
\mathfrak{P}_{j}
\right)
\right)_{j\in\bb{\mathbf{s}}}\right)&=
\sigma^{\mathbf{Pth}_{\boldsymbol{\mathcal{A}}}}\left(\left(
\mathrm{ip}^{(1,0)\sharp}_{s_{j}}\left(
\mathrm{sc}^{(0,1)}_{s_{j}}\left(
\mathfrak{P}_{j}
\right)\right)
\right)_{j\in\bb{\mathbf{s}}}\right)
\tag{1}
\\&=
\mathrm{ip}^{(1,0)\sharp}_{s}\left(
\sigma^{\mathbf{T}_{\Sigma}(X)}\left(
\left(
\mathrm{sc}^{(1,0)}_{s_{j}}\left(
\mathfrak{P}_{j}
\right)
\right)_{j\in\bb{\mathbf{s}}}
\right)
\right)
\tag{2}
\\&=
\mathrm{ip}^{(1,0)\sharp}_{s}\left(
\mathrm{sc}^{(1,0)}_{s}\left(
\sigma^{\mathbf{Pth}_{\boldsymbol{\mathcal{A}}}}\left(
\left(
\mathfrak{P}_{j}
\right)_{j\in\bb{\mathbf{s}}}
\right)
\right)
\right)
\tag{3}
\\&=
\mathrm{sc}^{0\mathbf{Pth}_{\boldsymbol{\mathcal{A}}}}_{s}\left(
\sigma^{\mathbf{Pth}_{\boldsymbol{\mathcal{A}}}}\left(
\left(
\mathfrak{P}_{s_{j}}
\right)_{j\in\bb{\mathbf{s}}}
\right)
\right).
\tag{4}
\end{align*}

In the just stated chain of equalities, the first equality unravels the description of the $0$-source according to Proposition~\ref{PPthCatAlg}; the second equality follows from Proposition~\ref{PSigma}; the third equality follows from Proposition~\ref{PHom}; finally the last equality recovers the description of the $0$-source according to Proposition~\ref{PPthCatAlg}.

Note that the desired equality follows from the fact that the $\mathrm{Ker}(\mathrm{CH}^{(1)})$-classes of the paths above are equal.

This completes the proof.
\end{proof}

\begin{proposition}\label{PCHVarA8} Let $(\mathbf{s},s)$ be an element in $S^{\star}\times S$, $\sigma$ an operation symbol in $\Sigma^{\boldsymbol{\mathcal{A}}}_{\mathbf{s},s}$ and $([\mathfrak{P}_{j}]_{s_{j}})_{j\in\bb{\mathbf{s}}}$, $([\mathfrak{Q}_{j}]_{s_{j}})_{j\in\bb{\mathbf{s}}}$ families of path classes in $[\mathrm{Pth}_{\boldsymbol{\mathcal{A}}}]_{\mathbf{s}}$, such that, for every $j\in\bb{\mathbf{s}}$, 
$$
\mathrm{sc}^{0[\mathbf{Pth}_{\boldsymbol{\mathcal{A}}}]}_{s_{j}}\left(
\left[
\mathfrak{Q}_{j}
\right]_{s_{j}}
\right)
=
\mathrm{tg}^{0[\mathbf{Pth}_{\boldsymbol{\mathcal{A}}}]}_{s_{j}}\left(
\left[
\mathfrak{P}_{j}
\right]_{s_{j}}
\right).
$$
Then the following equality holds
\allowdisplaybreaks
\begin{multline*}
\sigma^{[\mathbf{Pth}_{\boldsymbol{\mathcal{A}}}]}
\left(\left(
\left[\mathfrak{Q}_{j}\right]_{s_{j}}
\circ^{0[\mathbf{Pth}_{\boldsymbol{\mathcal{A}}}]}_{s_{j}}
\left[\mathfrak{P}_{j}\right]_{s_{j}}
\right)_{j\in\bb{\mathbf{s}}}
\right)
\\=
\sigma^{[\mathbf{Pth}_{\boldsymbol{\mathcal{A}}}]}
\left(\left(
\left[\mathfrak{Q}_{j}\right]_{s_{j}}
\right)_{j\in\bb{\mathbf{s}}}
\right)
\circ^{0[\mathbf{Pth}_{\boldsymbol{\mathcal{A}}}]}_{s}
\sigma^{[\mathbf{Pth}_{\boldsymbol{\mathcal{A}}}]}
\left(\left(
\left[\mathfrak{P}_{j}\right]_{s_{j}}
\right)_{j\in\bb{\mathbf{s}}}
\right).
\end{multline*}
\end{proposition}
\begin{proof}

We consider different possibilities according to the nature of $\sigma$. It could be the case that (1) $\mathbf{s}=\lambda$ or (2) $\mathbf{s}\in S^{\star}-\{\lambda\}$.

If~(1), i.e., if $\sigma$ is a constant operation symbol in $\Sigma_{\lambda,s}$ then the following chain of equalities hold
\allowdisplaybreaks
\begin{align*}
\sigma^{[\mathbf{Pth}_{\boldsymbol{\mathcal{A}}}]}
&=
\left[
\sigma^{\mathbf{Pth}_{\boldsymbol{\mathcal{A}}}}
\right]_{s}
\tag{1}
\\&=
\left[
\sigma^{\mathbf{Pth}_{\boldsymbol{\mathcal{A}}}}
\circ^{0\mathbf{Pth}_{\boldsymbol{\mathcal{A}}}}_{s}
\sigma^{\mathbf{Pth}_{\boldsymbol{\mathcal{A}}}}
\right]_{s}
\tag{2}
\\&=
\left[
\sigma^{\mathbf{Pth}_{\boldsymbol{\mathcal{A}}}}
\right]_{s}
\circ^{0[\mathbf{Pth}_{\boldsymbol{\mathcal{A}}}]}_{s}
\left[
\sigma^{\mathbf{Pth}_{\boldsymbol{\mathcal{A}}}}
\right]_{s}
\tag{3}
\\&=
\sigma^{[\mathbf{Pth}_{\boldsymbol{\mathcal{A}}}]}
\circ^{0[\mathbf{Pth}_{\boldsymbol{\mathcal{A}}}]}_{s}
\sigma^{[\mathbf{Pth}_{\boldsymbol{\mathcal{A}}}]}.
\tag{4}
\end{align*}

The first equality unravels the interpretation of the $\sigma$ operation in the many-sorted partial $\Sigma^{\boldsymbol{\mathcal{A}}}$-algebra $[\mathbf{Pth}_{\boldsymbol{\mathcal{A}}}]$ according to Proposition~\ref{PCHCatAlg}; the second equality follows from the fact that, according to Proposition~\ref{PPthCatAlg}, the interpretation of the $\sigma$ operation in the many-sorted partial $\Sigma^{\boldsymbol{\mathcal{A}}}$-algebra $\mathbf{Pth}_{\boldsymbol{\mathcal{A}}}$  is given by the $(1,0)$-identity path on $\sigma^{\mathbf{T}_{\Sigma}(X)}$, i.e., 
$
\sigma^{\mathbf{Pth}_{\boldsymbol{\mathcal{A}}}}=
\mathrm{ip}^{(1,0)\sharp}_{s}\left(\sigma^{\mathbf{T}_{\Sigma}(X)}\right)
$ and from the fact that, according to Proposition~\ref{PPthComp} the $(1,0)$-identity paths are idempotent for the $0$-composition. Note that this $0$-composition is well-defined by Proposition~\ref{PBasicEq}; the third equality recovers the $0$-composition operation in the many-sorted partial $\Sigma^{\boldsymbol{\mathcal{A}}}$-algebra $[\mathbf{Pth}_{\boldsymbol{\mathcal{A}}}]$ according to Proposition~\ref{PCHCatAlg}; finally, the last equality recovers the interpretation of the $\sigma$ operation in the many-sorted partial $\Sigma^{\boldsymbol{\mathcal{A}}}$-algebra $[\mathbf{Pth}_{\boldsymbol{\mathcal{A}}}]$ according to Proposition~\ref{PCHCatAlg}.

Case~(1) follows.

Regarding case~(2), i.e., if $\sigma\in\Sigma_{\mathbf{s},s}$ for some $\mathbf{s}\in S^{\star}-\{\lambda\}$, we consider the subcases~(2.1), where $(\mathfrak{P}_{j})_{j\in\bb{\mathbf{s}}}$ and $(\mathfrak{Q}_{j})_{j\in\bb{\mathbf{s}}}$ are families of $(1,0)$-identity paths; (2.2), where $(\mathfrak{P}_{j})_{j\in\bb{\mathbf{s}}}$ is a family of $(1,0)$-identity paths and $(\mathfrak{Q}_{j})_{j\in\bb{\mathbf{s}}}$ is not a family of $(1,0)$-identity paths;  (2.3), where $(\mathfrak{P}_{j})_{j\in\bb{\mathbf{s}}}$ is not a family of $(1,0)$-identity paths and $(\mathfrak{Q}_{j})_{j\in\bb{\mathbf{s}}}$ is a family of $(1,0)$-identity paths; and (2.4), where neither $(\mathfrak{P}_{j})_{j\in\bb{\mathbf{s}}}$ nor $(\mathfrak{Q}_{j})_{j\in\bb{\mathbf{s}}}$ are families of of $(1,0)$-identity paths.

If~(2.1), i.e., if $\sigma$ is a non-constant operation symbol in $\Sigma_{\mathbf{s},s}$ and $(\mathfrak{P}_{j})_{j\in\bb{\mathbf{s}}}$ and $(\mathfrak{Q}_{j})_{j\in\bb{\mathbf{s}}}$ are families of $(1,0)$-identity paths, we now proceed to describe the situation under consideration.

Let us recall that, for every $j\in\bb{\mathbf{s}}$, 
$$
\mathrm{sc}^{0[\mathbf{Pth}_{\boldsymbol{\mathcal{A}}}]}_{s_{j}}\left(
\left[
\mathfrak{Q}_{j}
\right]_{s_{j}}
\right)
=
\mathrm{tg}^{0[\mathbf{Pth}_{\boldsymbol{\mathcal{A}}}]}_{s_{j}}\left(
\left[
\mathfrak{P}_{j}
\right]_{s_{j}}
\right).
$$

Thus, following Proposition~\ref{PPthCatAlg} and Corollary~\ref{CCHUZId}, we conclude that, for every $j\in\bb{\mathbf{s}}$, there exists a term $P_{j}$ in $\mathrm{T}_{\Sigma}(X)_{s_{j}}$ for which
\begin{align*}
\mathfrak{P}_{j}=\mathrm{ip}^{(1,0)\sharp}_{s_{j}}(P_{j})=\mathfrak{Q}_{j}.
\tag{$\star$}
\end{align*}

The following chain of equalities holds
\begin{flushleft}
$\sigma^{\mathbf{Pth}_{\boldsymbol{\mathcal{A}}}}
\left(\left(
\mathfrak{Q}_{j}
\circ^{0\mathbf{Pth}_{\boldsymbol{\mathcal{A}}}}_{s_{j}}
\mathfrak{P}_{j}
\right)_{j\in\bb{\mathbf{s}}}
\right)$
\allowdisplaybreaks
\begin{align*}
\qquad
&=
\sigma^{\mathbf{Pth}_{\boldsymbol{\mathcal{A}}}}
\left(\left(
\mathrm{ip}^{(1,0)\sharp}_{s_{j}}\left(
P_{j}
\right)
\circ^{0\mathbf{Pth}_{\boldsymbol{\mathcal{A}}}}_{s_{j}}
\mathrm{ip}^{(1,0)\sharp}_{s_{j}}\left(
P_{j}
\right)
\right)_{j\in\bb{\mathbf{s}}}
\right)
\tag{1}
\\&=
\sigma^{\mathbf{Pth}_{\boldsymbol{\mathcal{A}}}}
\left(\left(
\mathrm{ip}^{(1,0)\sharp}_{s_{j}}\left(
P_{j}
\right)
\right)_{j\in\bb{\mathbf{s}}}
\right)
\tag{2}
\\&=
\mathrm{ip}^{(1,0)\sharp}_{s}\left(
\sigma^{\mathbf{T}_{\Sigma}(X)}
\left(\left(
P_{j}
\right)_{j\in\bb{\mathbf{s}}}
\right)
\right)
\tag{3}
\\&=
\mathrm{ip}^{(1,0)\sharp}_{s}\left(
\sigma^{\mathbf{T}_{\Sigma}(X)}
\left(\left(
P_{j}
\right)_{j\in\bb{\mathbf{s}}}
\right)
\right)
\circ_{s}^{0\mathbf{Pth}_{\boldsymbol{\mathcal{A}}}}
\mathrm{ip}^{(1,0)\sharp}_{s}\left(
\sigma^{\mathbf{T}_{\Sigma}(X)}
\left(\left(
P_{j}
\right)_{j\in\bb{\mathbf{s}}}
\right)
\right)
\tag{4}
\\&=
\sigma^{\mathbf{Pth}_{\boldsymbol{\mathcal{A}}}}
\left(\left(
\mathrm{ip}^{(1,0)\sharp}_{s_{j}}\left(
P_{j}
\right)
\right)_{j\in\bb{\mathbf{s}}}
\right)
\circ_{s}^{0\mathbf{Pth}_{\boldsymbol{\mathcal{A}}}}
\sigma^{\mathbf{Pth}_{\boldsymbol{\mathcal{A}}}}
\left(\left(
\mathrm{ip}^{(1,0)\sharp}_{s_{j}}\left(
P_{j}
\right)
\right)_{j\in\bb{\mathbf{s}}}
\right)
\tag{5}
\\&=
\sigma^{\mathbf{Pth}_{\boldsymbol{\mathcal{A}}}}
\left(\left(
\mathfrak{Q}_{j}
\right)_{j\in\bb{\mathbf{s}}}
\right)
\circ_{s}^{0\mathbf{Pth}_{\boldsymbol{\mathcal{A}}}}
\sigma^{\mathbf{Pth}_{\boldsymbol{\mathcal{A}}}}
\left(\left(
\mathfrak{P}_{j}
\right)_{j\in\bb{\mathbf{s}}}
\right).
\tag{6}
\end{align*}
\end{flushleft}

The first equality applies the equality presented in $(\star)$; the second equality follows from the fact that $(1,0)$-identity paths are idempotent for the $0$-composition; the third equality follows from the fact that, according to Proposition~\ref{PIpHom}, $\mathrm{ip}^{(1,0)\sharp}$ is a $\Sigma$-homomorphism; the fourth equality follows from the fact that $(1,0)$-identity paths are idempotent for the $0$-composition; the fifth equality follows from the fact that, according to Proposition~\ref{PIpHom}, $\mathrm{ip}^{(1,0)\sharp}$ is a $\Sigma$-homomorphism; finally, the last equality recovers the description, for every $j\in\bb{\mathbf{s}}$, of the paths $\mathfrak{P}_{j}$ and $\mathfrak{Q}_{j}$, according to the equality presented in $(\star)$.

Consequently, the following chain of equalities holds
\begin{flushleft}
$\sigma^{[\mathbf{Pth}_{\boldsymbol{\mathcal{A}}}]}
\left(\left(
\left[\mathfrak{Q}_{j}\right]_{s_{j}}
\circ^{0[\mathbf{Pth}_{\boldsymbol{\mathcal{A}}}]}_{s_{j}}
\left[\mathfrak{P}_{j}\right]_{s_{j}}
\right)_{j\in\bb{\mathbf{s}}}
\right)$
\allowdisplaybreaks
\begin{align*}
\qquad
&=
\left[
\sigma^{\mathbf{Pth}_{\boldsymbol{\mathcal{A}}}}
\left(\left(
\mathfrak{Q}_{j}
\circ^{0\mathbf{Pth}_{\boldsymbol{\mathcal{A}}}}_{s_{j}}
\mathfrak{P}_{j}
\right)_{j\in\bb{\mathbf{s}}}
\right)
\right]_{s}
\tag{1}
\\&=
\left[
\sigma^{\mathbf{Pth}_{\boldsymbol{\mathcal{A}}}}
\left(\left(
\mathfrak{Q}_{j}
\right)_{j\in\bb{\mathbf{s}}}
\right)
\circ_{s}^{0\mathbf{Pth}_{\boldsymbol{\mathcal{A}}}}
\sigma^{\mathbf{Pth}_{\boldsymbol{\mathcal{A}}}}
\left(\left(
\mathfrak{P}_{j}
\right)_{j\in\bb{\mathbf{s}}}
\right)
\right]_{s}
\tag{2}
\\&=
\sigma^{[\mathbf{Pth}_{\boldsymbol{\mathcal{A}}}]}
\left(\left(
\left[\mathfrak{Q}_{j}\right]_{s_{j}}
\right)_{j\in\bb{\mathbf{s}}}
\right)
\circ^{0[\mathbf{Pth}_{\boldsymbol{\mathcal{A}}}]}_{s}
\sigma^{[\mathbf{Pth}_{\boldsymbol{\mathcal{A}}}]}
\left(\left(
\left[\mathfrak{P}_{j}\right]_{s_{j}}
\right)_{j\in\bb{\mathbf{s}}}
\right).
\tag{3}
\end{align*}
\end{flushleft}

The first equality unravels the interpretation of the operations $\sigma$ and $\circ^{0}_{s}$ in the many-sorted partial $\Sigma^{\boldsymbol{\mathcal{A}}}$-algebra $[\mathbf{Pth}_{\boldsymbol{\mathcal{A}}}]$, according to Proposition~\ref{PCHCatAlg}; the second equality follows from the previous chain of equalities; finally, the last equality recovers the interpretation of the operations $\sigma$ and $\circ^{0}_{s}$ in the many-sorted partial $\Sigma^{\boldsymbol{\mathcal{A}}}$-algebra $[\mathbf{Pth}_{\boldsymbol{\mathcal{A}}}]$, according to Proposition~\ref{PCHCatAlg}.

Therefore, case~(2.1) follows.

If~(2.2), i.e., if $\sigma$ is a non-constant operation symbol in $\Sigma_{\mathbf{s},s}$ and $(\mathfrak{P}_{j})_{j\in\bb{\mathbf{s}}}$ is not a family of $(1,0)$-identity paths and $(\mathfrak{Q}_{j})_{j\in\bb{\mathbf{s}}}$ is a family of $(1,0)$-identity paths, we now proceed to describe the situation under consideration.

Let us recall that for every $j\in\bb{\mathbf{s}}$,
$$
\mathrm{sc}^{0[\mathbf{Pth}_{\boldsymbol{\mathcal{A}}}]}_{s_{j}}\left(
\left[
\mathfrak{Q}_{j}
\right]_{s_{j}}
\right)
=
\mathrm{tg}^{0[\mathbf{Pth}_{\boldsymbol{\mathcal{A}}}]}_{s_{j}}\left(
\left[
\mathfrak{P}_{j}
\right]_{s_{j}}
\right).
$$

Therefore, taking into account Proposition~\ref{PPthCatAlg}, Corollary~\ref{CCHUZId} and the fact  that, by assumption, $(\mathfrak{Q}_{j})_{j\in\bb{\mathbf{s}}}$ is a family of $(1,0)$-identity paths, we conclude that, for every $j\in\bb{\mathbf{s}}$, the following equality holds
\begin{align*}
\mathfrak{Q}_{j}=
\mathrm{ip}^{(1,0)\sharp}_{s_{j}}\left(
\mathrm{tg}^{(0,1)}_{s_{j}}\left(
\mathfrak{P}_{j}
\right)
\right)
\tag{A}
\end{align*}

Thus, following Proposition~\ref{PPthComp}, for every $j\in\bb{\mathbf{s}}$, the following equality holds
\begin{align*}
\mathfrak{Q}_{j}
\circ^{0\mathbf{Pth}_{\boldsymbol{\mathcal{A}}}}_{s_{j}}
\mathfrak{P}_{j}
=
\mathfrak{P}_{j}.
\tag{B}
\end{align*}

Moreover, the following chain of equalities holds
\allowdisplaybreaks
\begin{align*}
\sigma^{\mathbf{Pth}_{\boldsymbol{\mathcal{A}}}}
\left(\left(
\mathfrak{Q}_{j}
\right)_{j\in\bb{\mathbf{s}}}
\right)&=
\sigma^{\mathbf{Pth}_{\boldsymbol{\mathcal{A}}}}
\left(\left(
\mathrm{ip}^{(1,0)\sharp}_{s_{j}}\left(
\mathrm{tg}^{(0,1)}_{s_{j}}\left(
\mathfrak{P}_{j}
\right)
\right)
\right)_{j\in\bb{\mathbf{s}}}
\right)
\tag{1}
\\&=
\mathrm{ip}^{(1,0)\sharp}_{s}\left(
\sigma^{\mathbf{T}_{\Sigma}(X)}
\left(\left(
\mathrm{tg}^{(0,1)}_{s_{j}}\left(
\mathfrak{P}_{j}
\right)
\right)_{j\in\bb{\mathbf{s}}}
\right)
\right)
\tag{2}
\\&=
\mathrm{ip}^{(1,0)\sharp}_{s}\left(
\mathrm{tg}^{(0,1)}_{s}\left(
\sigma^{\mathbf{Pth}_{\boldsymbol{\mathcal{A}}}}
\left(\left(
\mathfrak{P}_{j}
\right)_{j\in\bb{\mathbf{s}}}
\right)
\right)
\right).
\tag{3}
\end{align*}

The first equality follows from the family of equalities in (A); the second equality follows from the fact that $\mathrm{ip}^{(1,0)\sharp}$ is a $\Sigma$-homomorphism, according to Proposition~\ref{PIpHom}; finally, the last equality follows from the fact that $\mathrm{tg}^{(0,1)}$ is a $\Sigma$-homomorphism, according to Proposition~\ref{PHom}.

Thus, following Proposition~\ref{PPthComp}, the following equality holds
\begin{align*}
\sigma^{\mathbf{Pth}_{\boldsymbol{\mathcal{A}}}}
\left(\left(
\mathfrak{Q}_{j}
\right)_{j\in\bb{\mathbf{s}}}
\right)
\circ^{0\mathbf{Pth}_{\boldsymbol{\mathcal{A}}}}_{s}
\sigma^{\mathbf{Pth}_{\boldsymbol{\mathcal{A}}}}
\left(\left(
\mathfrak{P}_{j}
\right)_{j\in\bb{\mathbf{s}}}
\right)
=
\sigma^{\mathbf{Pth}_{\boldsymbol{\mathcal{A}}}}
\left(\left(
\mathfrak{P}_{j}
\right)_{j\in\bb{\mathbf{s}}}
\right).
\tag{C}
\end{align*}

Consequently, the following chain of equalities holds
\begin{flushleft}
$\sigma^{[\mathbf{Pth}_{\boldsymbol{\mathcal{A}}}]}
\left(\left(
\left[\mathfrak{Q}_{j}\right]_{s_{j}}
\circ^{0[\mathbf{Pth}_{\boldsymbol{\mathcal{A}}}]}_{s_{j}}
\left[\mathfrak{P}_{j}\right]_{s_{j}}
\right)_{j\in\bb{\mathbf{s}}}
\right)$
\allowdisplaybreaks
\begin{align*}
\qquad
&=
\sigma^{[\mathbf{Pth}_{\boldsymbol{\mathcal{A}}}]}
\left(\left(
\left[
\mathfrak{Q}_{j}
\circ^{0\mathbf{Pth}_{\boldsymbol{\mathcal{A}}}}_{s_{j}}
\mathfrak{P}_{j}
\right]_{s_{j}}
\right)_{j\in\bb{\mathbf{s}}}
\right)
\tag{1}
\\&=
\sigma^{[\mathbf{Pth}_{\boldsymbol{\mathcal{A}}}]}
\left(\left(
\left[
\mathfrak{P}_{j}
\right]_{s_{j}}
\right)_{j\in\bb{\mathbf{s}}}
\right)
\tag{2}
\\&=
\left[
\sigma^{\mathbf{Pth}_{\boldsymbol{\mathcal{A}}}}
\left(\left(
\mathfrak{P}_{j}
\right)_{j\in\bb{\mathbf{s}}}
\right)
\right]_{s}
\tag{3}
\\&=
\left[
\sigma^{\mathbf{Pth}_{\boldsymbol{\mathcal{A}}}}
\left(\left(
\mathfrak{Q}_{j}
\right)_{j\in\bb{\mathbf{s}}}
\right)
\circ^{0\mathbf{Pth}_{\boldsymbol{\mathcal{A}}}}_{s}
\sigma^{\mathbf{Pth}_{\boldsymbol{\mathcal{A}}}}
\left(\left(
\mathfrak{P}_{j}
\right)_{j\in\bb{\mathbf{s}}}
\right)
\right]_{s}
\tag{4}
\\&=
\sigma^{[\mathbf{Pth}_{\boldsymbol{\mathcal{A}}}]}
\left(\left(
\left[\mathfrak{Q}_{j}\right]_{s_{j}}
\right)_{j\in\bb{\mathbf{s}}}
\right)
\circ^{0[\mathbf{Pth}_{\boldsymbol{\mathcal{A}}}]}_{s}
\sigma^{[\mathbf{Pth}_{\boldsymbol{\mathcal{A}}}]}
\left(\left(
\left[\mathfrak{P}_{j}\right]_{s_{j}}
\right)_{j\in\bb{\mathbf{s}}}
\right).
\tag{5}
\end{align*}
\end{flushleft}

The first equality unravels the interpretation of the operation $\circ^{0}_{s}$ in the many-sorted partial $\Sigma^{\boldsymbol{\mathcal{A}}}$-algebra $[\mathbf{Pth}_{\boldsymbol{\mathcal{A}}}]$, according to Proposition~\ref{PCHCatAlg};  the second equality  follows from the family of equalities in (B); the third equality unravels the interpretation of the operation $\sigma$ in the many-sorted partial $\Sigma^{\boldsymbol{\mathcal{A}}}$-algebra $[\mathbf{Pth}_{\boldsymbol{\mathcal{A}}}]$, according to Proposition~\ref{PCHCatAlg}; the fourth equality follows from equality (C); finally, the last equality recovers the interpretation of the oeprations $\sigma$ and $\circ^{0}_{s}$ in the many-sorted partial $\Sigma^{\boldsymbol{\mathcal{A}}}$-algebra $[\mathbf{Pth}_{\boldsymbol{\mathcal{A}}}]$, according to Proposition~\ref{PCHCatAlg}.

Therefore, case~(2,2) follows.

If~(2.3), i.e., if $\sigma$ is a non-constant operation symbol in $\Sigma_{\mathbf{s},s}$ and $(\mathfrak{P}_{j})_{j\in\bb{\mathbf{s}}}$ is a family of $(1,0)$-identity paths and $(\mathfrak{Q}_{j})_{j\in\bb{\mathbf{s}}}$ is not a family of $(1,0)$-identity paths, we  proceed as in case~(2.2).

If~(2.4), i.e., if $\sigma$ is a non-constant operation symbol in $\Sigma_{\mathbf{s},s}$ and neither $(\mathfrak{P}_{j})_{j\in\bb{\mathbf{s}}}$ nor  $(\mathfrak{Q}_{j})_{j\in\bb{\mathbf{s}}}$ are families of $(1,0)$-identity paths, we now proceed to describe the situation under consideration.

According to Corollary~\ref{CPthWB}, we have that
\begin{enumerate}
\item[(i)] $\sigma^{\mathbf{Pth}_{\boldsymbol{\mathcal{A}}}}
((
\mathfrak{P}_{j}
)_{j\in\bb{\mathbf{s}}}
)
$ is an echelonless path. Moreover, according to Proposition~\ref{PRecov}, the path extraction algorithm from Lemma~\ref{LPthExtract} applied to the path $\sigma^{\mathbf{Pth}_{\boldsymbol{\mathcal{A}}}}
((
\mathfrak{P}_{j}
)_{j\in\bb{\mathbf{s}}}
)
$ retrieves the original family $(
\mathfrak{P}_{j}
)_{j\in\bb{\mathbf{s}}}
$
;
\item[(ii)] $\sigma^{\mathbf{Pth}_{\boldsymbol{\mathcal{A}}}}
((
\mathfrak{Q}_{j}
)_{j\in\bb{\mathbf{s}}}
)
$ is an echelonless path. Moreover, according to Proposition~\ref{PRecov}, the path extraction algorithm from Lemma~\ref{LPthExtract} applied to the path $\sigma^{\mathbf{Pth}_{\boldsymbol{\mathcal{A}}}}
((
\mathfrak{Q}_{j}
)_{j\in\bb{\mathbf{s}}}
)
$ retrieves the original family $(
\mathfrak{Q}_{j}
)_{j\in\bb{\mathbf{s}}}
$.
\end{enumerate}

Since we are assuming that,  for every $j\in\bb{\mathbf{s}}$, the following equality holds
$$
\mathrm{sc}^{0[\mathbf{Pth}_{\boldsymbol{\mathcal{A}}}]}_{s_{j}}\left(
\left[
\mathfrak{Q}_{j}
\right]_{s_{j}}
\right)
=
\mathrm{tg}^{0[\mathbf{Pth}_{\boldsymbol{\mathcal{A}}}]}_{s_{j}}\left(
\left[
\mathfrak{P}_{j}
\right]_{s_{j}}
\right),
$$
we conclude, on one hand, that the family of $0$-compositions $(\mathfrak{Q}_{j}\circ^{0\mathbf{Pth}_{\boldsymbol{\mathcal{A}}}}_{s_{j}}\mathfrak{P}_{j})_{j\in\bb{\mathbf{s}}}$ is well-defined. Moreover, taking account the assumptions of case~(2.4), this is a family not entirely compose of $(1,0)$-identity paths.  According to Corollary~\ref{CPthWB}, we have that
\begin{enumerate}
\item[(iii)] $\sigma^{\mathbf{Pth}_{\boldsymbol{\mathcal{A}}}}
((
\mathfrak{Q}_{j}\circ^{0\mathbf{Pth}_{\boldsymbol{\mathcal{A}}}}_{s_{j}}\mathfrak{P}_{j}
)_{j\in\bb{\mathbf{s}}}
)
$ is an echelonless path. Moreover, according to Proposition~\ref{PRecov}, the path extraction algorithm from Lemma~\ref{LPthExtract} applied to the path $\sigma^{\mathbf{Pth}_{\boldsymbol{\mathcal{A}}}}
((
\mathfrak{Q}_{j}\circ^{0\mathbf{Pth}_{\boldsymbol{\mathcal{A}}}}_{s_{j}}\mathfrak{P}_{j}
)_{j\in\bb{\mathbf{s}}}
)
$ retrieves the original family $(\mathfrak{Q}_{j}\circ^{0\mathbf{Pth}_{\boldsymbol{\mathcal{A}}}}_{s_{j}}\mathfrak{P}_{j})_{j\in\bb{\mathbf{s}}}$
.
\end{enumerate}

On the other hand,  we conclude in virtue of Propositions~\ref{PPthCatAlg},~\ref{PHom} and~\ref{PPthCatAlg} that the path given by the $0$-composition
$$
\sigma^{\mathbf{Pth}_{\boldsymbol{\mathcal{A}}}}
\left(\left(
\mathfrak{Q}_{j}
\right)_{j\in\bb{\mathbf{s}}}
\right)
\circ^{0\mathbf{Pth}_{\boldsymbol{\mathcal{A}}}}_{s}
\sigma^{\mathbf{Pth}_{\boldsymbol{\mathcal{A}}}}
\left(\left(
\mathfrak{P}_{j}
\right)_{j\in\bb{\mathbf{s}}}
\right)
$$ 
is well-defined.

Taking into account items (i) and (ii) above, we conclude that
\begin{enumerate}
\item[(iv)] $
\sigma^{\mathbf{Pth}_{\boldsymbol{\mathcal{A}}}}
((
\mathfrak{Q}_{j}
)_{j\in\bb{\mathbf{s}}}
)
\circ^{0\mathbf{Pth}_{\boldsymbol{\mathcal{A}}}}_{s}
\sigma^{\mathbf{Pth}_{\boldsymbol{\mathcal{A}}}}
((
\mathfrak{P}_{j}
)_{j\in\bb{\mathbf{s}}}
)
$ is an echelonless path.
\end{enumerate}

If we apply the path extraction algorithm from Lemma~\ref{LPthExtract} to the echelonless path $\sigma^{\mathbf{Pth}_{\boldsymbol{\mathcal{A}}}}
((
\mathfrak{Q}_{j}
)_{j\in\bb{\mathbf{s}}}
)
\circ^{0\mathbf{Pth}_{\boldsymbol{\mathcal{A}}}}_{s}
\sigma^{\mathbf{Pth}_{\boldsymbol{\mathcal{A}}}}
((
\mathfrak{P}_{j}
)_{j\in\bb{\mathbf{s}}}
)$, we have that, for every $j\in\bb{\mathbf{s}}$, the $j$-th component we can extract from it is given by
\begin{flushleft}
$\left(
\sigma^{\mathbf{Pth}_{\boldsymbol{\mathcal{A}}}}
\left(\left(
\mathfrak{Q}_{j}
\right)_{j\in\bb{\mathbf{s}}}
\right)
\circ^{0\mathbf{Pth}_{\boldsymbol{\mathcal{A}}}}_{s}
\sigma^{\mathbf{Pth}_{\boldsymbol{\mathcal{A}}}}
\left(\left(
\mathfrak{P}_{j}
\right)_{j\in\bb{\mathbf{s}}}
\right)
\right)_{j}$
\allowdisplaybreaks
\begin{align*}
\qquad
&=
\left(
\sigma^{\mathbf{Pth}_{\boldsymbol{\mathcal{A}}}}
\left(\left(
\mathfrak{Q}_{j}
\right)_{j\in\bb{\mathbf{s}}}
\right)
\right)_{j}
\circ^{0\mathbf{Pth}_{\boldsymbol{\mathcal{A}}}}_{s_{j}}
\left(
\sigma^{\mathbf{Pth}_{\boldsymbol{\mathcal{A}}}}
\left(\left(
\mathfrak{P}_{j}
\right)_{j\in\bb{\mathbf{s}}}
\right)
\right)_{j}
\tag{1}
\\&=
\mathfrak{Q}_{j}
\circ^{0\mathbf{Pth}_{\boldsymbol{\mathcal{A}}}}_{s_{j}}
\mathfrak{P}_{j}.
\tag{2}
\end{align*}
\end{flushleft}

In the just stated chain of equalities, the first equality follows from the fact that the $j$-th component we can extract from the $0$-composition is given by the $0$-composition of the corresponding $j$-th components; finally, the last equality follows from items (i) and (ii) above.

It follows that the families of paths we can extract from the different echelonless paths considered in items~(iii) and~(iv) are equal and they are given by
$$
\left(
\mathfrak{Q}_{j}\circ^{0\mathbf{Pth}_{\boldsymbol{\mathcal{A}}}}_{s_{j}}\mathfrak{P}_{j}
\right)_{j\in\bb{\mathbf{s}}}.
$$

As a consequence, the echelonless paths considered in items~(iii) and~(iv) have the same value for the the Curry-Howard mapping. 
\begin{flushleft}
$
\mathrm{CH}^{(1)}_{s}\left(
\sigma^{\mathbf{Pth}_{\boldsymbol{\mathcal{A}}}}
\left(\left(
\mathfrak{Q}_{j}
\right)_{j\in\bb{\mathbf{s}}}
\right)
\circ^{0\mathbf{Pth}_{\boldsymbol{\mathcal{A}}}}_{s}
\sigma^{\mathbf{Pth}_{\boldsymbol{\mathcal{A}}}}
\left(\left(
\mathfrak{P}_{j}
\right)_{j\in\bb{\mathbf{s}}}
\right)
\right)
$
\allowdisplaybreaks
\begin{align*}
\qquad&=
\sigma^{\mathbf{T}_{\Sigma^{\boldsymbol{\mathcal{A}}}}(X)}
\left(\left(
\mathrm{CH}^{(1)}_{s_{j}}
\left(
\mathfrak{Q}_{j}
\circ^{0\mathbf{Pth}_{\boldsymbol{\mathcal{A}}}}_{s_{j}}
\mathfrak{P}_{j}
\right)
\right)_{j\in\bb{\mathbf{s}}}
\right)
\tag{1}
\\&=
\mathrm{CH}^{(1)}_{s}\left(
\sigma^{\mathbf{Pth}_{\boldsymbol{\mathcal{A}}}}
\left(\left(
\mathfrak{Q}_{j}\circ^{0\mathbf{Pth}_{\boldsymbol{\mathcal{A}}}}_{s_{j}}\mathfrak{P}_{j}
\right)_{j\in\bb{\mathbf{s}}}
\right)
\right).
\tag{2}
\end{align*}
\end{flushleft}

The first equality follows from Definition~\ref{DCH} and item~(iv); the second equality follows from Definition~\ref{DCH} and item (iii).

Thus, the following chain of equalities holds
\begin{flushleft}
$\sigma^{[\mathbf{Pth}_{\boldsymbol{\mathcal{A}}}]}
\left(\left(
\left[\mathfrak{Q}_{j}\right]_{s_{j}}
\circ^{0[\mathbf{Pth}_{\boldsymbol{\mathcal{A}}}]}_{s_{j}}
\left[\mathfrak{P}_{j}\right]_{s_{j}}
\right)_{j\in\bb{\mathbf{s}}}
\right)$
\allowdisplaybreaks
\begin{align*}
\qquad&=
\left[
\sigma^{\mathbf{Pth}_{\boldsymbol{\mathcal{A}}}}
\left(\left(
\mathfrak{Q}_{j}\circ^{0\mathbf{Pth}_{\boldsymbol{\mathcal{A}}}}_{s_{j}}\mathfrak{P}_{j}
\right)_{j\in\bb{\mathbf{s}}}
\right)
\right]_{s}
\tag{1}
\\&=
\left[
\sigma^{\mathbf{Pth}_{\boldsymbol{\mathcal{A}}}}
\left(\left(
\mathfrak{Q}_{j}
\right)_{j\in\bb{\mathbf{s}}}
\right)
\circ^{0\mathbf{Pth}_{\boldsymbol{\mathcal{A}}}}_{s}
\sigma^{\mathbf{Pth}_{\boldsymbol{\mathcal{A}}}}
\left(\left(
\mathfrak{P}_{j}
\right)_{j\in\bb{\mathbf{s}}}
\right)
\right]_{s}
\tag{2}
\\&=
\sigma^{[\mathbf{Pth}_{\boldsymbol{\mathcal{A}}}]}
\left(\left(
\left[\mathfrak{Q}_{j}\right]_{s_{j}}
\right)_{j\in\bb{\mathbf{s}}}
\right)
\circ^{0[\mathbf{Pth}_{\boldsymbol{\mathcal{A}}}]}_{s}
\sigma^{[\mathbf{Pth}_{\boldsymbol{\mathcal{A}}}]}
\left(\left(
\left[\mathfrak{P}_{j}\right]_{s_{j}}
\right)_{j\in\bb{\mathbf{s}}}
\right).
\tag{3}
\end{align*}
\end{flushleft}

The first equality unravels the interpretation of the operations $\sigma$ and $\circ^{0}_{s}$ in the many-sorted partial $\Sigma^{\boldsymbol{\mathcal{A}}}$-algebra $[\mathbf{Pth}_{\boldsymbol{\mathcal{A}}}]$, according to Proposition~\ref{PCHCatAlg}; the second equality follows from the previous discussion; finally, the last equality recovers the interpretation of the operations $\sigma$ and $\circ^{0}_{s}$ in the many-sorted partial $\Sigma^{\boldsymbol{\mathcal{A}}}$-algebra $[\mathbf{Pth}_{\boldsymbol{\mathcal{A}}}]$, according to Proposition~\ref{PCHCatAlg}.

Therefore, case~(2.4) follows.

This completes the proof.
\end{proof}

\begin{restatable}{proposition}{PCHCtyAlg}
\label{PCHCtyAlg} $[\mathsf{Pth}_{\boldsymbol{\mathcal{A}}}]$ is an $S$-sorted categorial $\Sigma$-algebra.
\end{restatable}
\begin{proof}
That $[\mathsf{Pth}_{\boldsymbol{\mathcal{A}}}]$ is an $S$-sorted category was already proven in Proposition~\ref{PCHCat}.  We are only left to prove that, for every $(\mathbf{s},s)$ in $S^{\star}\times S$, $\sigma^{[\boldsymbol{\mathsf{Pth}}_{\boldsymbol{\mathcal{A}}}]}$ is a functor from the product category $[\mathsf{Pth}_{\boldsymbol{\mathcal{A}}}]_{\mathbf{s}}$ to the category $[\mathsf{Pth}_{\boldsymbol{\mathcal{A}}}]_{s}$. But this  follows from Propositions~\ref{PCHVarA7} and~\ref{PCHVarA8}.

This completes the proof.
\end{proof}

\section{
\texorpdfstring
{An Artinian order on $\coprod[\mathrm{Pth}_{\boldsymbol{\mathcal{A}}}]$}
{An Artinian order on the quotient of paths}
}

In this section we define on $\coprod[\mathrm{Pth}_{\boldsymbol{\mathcal{A}}}]$, the coproduct of $[\mathrm{Pth}_{\boldsymbol{\mathcal{A}}}]$---formed by all labelled path classes $([\mathfrak{P}]_{s},s)$ with $s\in S$ and 
$[\mathfrak{P}]_{s}\in[\mathrm{Pth}_{\boldsymbol{\mathcal{A}}}]_{s}$---, an Artinian order (see Definition~\ref{DPosArt}) $\leq_{[\mathbf{Pth}_{\boldsymbol{\mathcal{A}}}]}$, which will allow us to justify both proofs by Artinian induction and definitions by Artinian recursion, which we will use in subsequent chapters. Moreover, we prove that several mappings from and to $\coprod[\mathrm{Pth}_{\boldsymbol{\mathcal{A}}}]$ are 
order-preserving, order-reflecting or order embeddings.

\begin{restatable}{definition}{DCHOrd}
\label{DCHOrd} 
\index{partial order!first-order!$\leq_{[\mathbf{Pth}_{\boldsymbol{\mathcal{A}}}]}$}
Let $\leq_{[\mathbf{Pth}_{\boldsymbol{\mathcal{A}}}]}$ be the binary relation defined on $\coprod[\mathrm{Pth}_{\boldsymbol{\mathcal{A}}}]$ containing every pair
$(([\mathfrak{Q}]_{t},t),([\mathfrak{P}]_{s},s))$ in $(\coprod[\mathrm{Pth}_{\boldsymbol{\mathcal{A}}}])^{2}$ satisfying that
$$
\exists\, \mathfrak{Q}'\in \left[
\mathfrak{Q}
\right]_{t},\,\, 
\exists\, \mathfrak{P}'\in 
\left[
\mathfrak{P}
\right]_{s}\, 
\left(
\left(
\mathfrak{Q}',t
\right)
\leq_{\mathbf{Pth}_{\boldsymbol{\mathcal{A}}}}
\left(
\mathfrak{P}',s
\right)
\right).
$$

That is, $([\mathfrak{Q}]_{t},t)$ $\leq_{[\mathbf{Pth}_{\boldsymbol{\mathcal{A}}}]}$-precedes $([\mathfrak{P}]_{s},s)$ if there exists a pair of representative paths $\mathfrak{Q}'\in [\mathfrak{Q}]_{t}$ and $\mathfrak{P}'\in [\mathfrak{P}]_{s}$ for which $(\mathfrak{Q}',t)$ $\leq_{\mathbf{Pth}_{\boldsymbol{\mathcal{A}}}}$-precedes $(\mathfrak{P}',s)$.
\end{restatable}

The aim of this section is to prove that $\leq_{[\mathbf{Pth}_{\boldsymbol{\mathcal{A}}}]}$ is an Artinian  order on $\coprod[\mathrm{Pth}_{\boldsymbol{\mathcal{A}}}]$. To prove it we need to state several lemmas, which, ultimately, are based on the work by Ganesamoorthy and Karpagaval in~\cite{GK15} (let us point out that, for quotients of preordered sets, Bourbaki, in~\cite{nB70}, Exercise 2, \S 1, Chapter III, considers a property of an equivalence relation on a preordered set, that of being ``weakly compatible'', which is the same as one of those considered by Ganesamoorthy and Karpagaval in~\cite{GK15}; another of the conditions considered by these authors being the dual of Bourbaki's condition. Moreover, for quotients of ordered sets Blyth, in~\cite{bl05}, on p. 41, considers a property of an equivalence relation on an ordered set, the  ``link property'', which is the same as one of those considered in~\cite{GK15}).

In the following lemma we show that, if we are given two paths in $\mathrm{Ker}(\mathrm{CH}^{(1)})$ and one of these paths is smaller than another path with respect to the order $\leq_{\mathbf{Pth}_{\boldsymbol{\mathcal{A}}}}$, then we can complete this inequality on the other side with an equivalent path.

\begin{lemma}\label{LCHOrdI} Let $t,s$ be sorts in $S$, let $\mathfrak{Q}',\mathfrak{Q}$ be paths in $\mathrm{Pth}_{\boldsymbol{\mathcal{A}},t}$ and let $\mathfrak{P}$ be a path in $\mathrm{Pth}_{\boldsymbol{\mathcal{A}},s}$. If  $(\mathfrak{Q},t)\leq_{\mathbf{Pth}_{\boldsymbol{\mathcal{A}}}}(\mathfrak{P},s)$ and $\mathfrak{Q}'\in[\mathfrak{Q}]_{t}$ then there exists a path $\mathfrak{P}'\in [\mathfrak{P}]_{s}$ such  that  
$$\left(
\mathfrak{Q}',t
\right)\leq_{\mathbf{Pth}_{\boldsymbol{\mathcal{A}}}}
\left(
\mathfrak{P}',s
\right).$$
\end{lemma}
\begin{proof}
Let us recall from Remark~\ref{ROrd} that $(\mathfrak{Q},t)\leq_{\mathbf{Pth}_{\boldsymbol{\mathcal{A}}}}(\mathfrak{P},s)$ if and only if $s=t$ and $\mathfrak{Q}=\mathfrak{P}$ or there exists a natural number $m\in\mathbb{N}-\{0\}$, a word $\mathbf{w}\in S^{\star}$ of length $\bb{\mathbf{w}}=m+1$, and a family of paths $(\mathfrak{R}_{k})$ in $\mathrm{Pth}_{\boldsymbol{\mathcal{A}},\mathbf{w}}$ such that $w_{0}=t$, $\mathfrak{R}_{0}=\mathfrak{Q}$, $w_{m}=s$, $\mathfrak{R}_{m}=\mathfrak{P}$ and for every $k\in m$, $(\mathfrak{R}_{k}, w_{k})\prec_{\mathbf{Pth}_{\boldsymbol{\mathcal{A}}}} (\mathfrak{R}_{k+1}, w_{k+1})$.

The lemma holds trivially in case $s=t$ and $\mathfrak{Q}=\mathfrak{P}$. Therefore, it remains to prove the other case. 

We will prove the lemma by induction on $m\in\mathbb{N}-\{0\}$.

\textsf{Base step of the induction.}

For $m=1$ we have that $(\mathfrak{Q},t)\prec_{\mathbf{Pth}_{\boldsymbol{\mathcal{A}}}}(\mathfrak{P},s)$, hence following Definition~\ref{DOrd} we are in one of the following situations
\begin{enumerate}
\item $\mathfrak{P}$ and $\mathfrak{Q}$ are $(1,0)$-identity paths of the form
\begin{align*}
\mathfrak{P}&=
\mathrm{ip}^{(1,0)\sharp}_{s}(P),
&
\mathfrak{Q}&=
\mathrm{ip}^{(1,0)\sharp}_{t}(Q),
\end{align*}
for some terms $P\in\mathrm{T}_{\Sigma}(X)_{s}$ and $Q\in\mathrm{T}_{\Sigma}(X)_{t}$ and the following inequality holds
$$
\left(
Q, t
\right)
<_{\mathbf{T}_{\Sigma}(X)}
\left(
P,s
\right),
$$
where $\leq_{\mathbf{T}_{\Sigma}(X)}$ is the Artinian partial order on $\coprod\mathrm{T}_{\Sigma}(X)$ introduced in Remark~\ref{RTermOrd}.

\item $\mathfrak{P}$ is a path of length strictly greater than one  containing at least one echelon, and if its first echelon occurs at position $i\in\bb{\mathfrak{P}}$, then
\subitem if $i=0$, then $\mathfrak{Q}$ is equal to $\mathfrak{P}^{0,0}$ or $\mathfrak{P}^{1,\bb{\mathfrak{P}}-1}$,
\subitem if $i>$, then $\mathfrak{Q}$ is equal to $\mathfrak{P}^{0,i-1}$ or $\mathfrak{P}^{i,\bb{\mathfrak{P}}-1}$;
\item $\mathfrak{P}$ is an echelonless path and $\mathfrak{Q}$ is one of the paths we can extract from $\mathfrak{P}$ in virtue of Lemma~\ref{LPthExtract}.
\end{enumerate}

If~(1), we note that both $\mathfrak{P}$ and $\mathfrak{Q}$ are $(1,0)$-identity paths. Then, if $\mathfrak{Q}'\in[\mathfrak{Q}]_{t}$ then we conclude, in virtue of Corollary~\ref{CCHUZId}, that $\mathfrak{Q}'=\mathfrak{Q}$. This case follows easily since $\mathfrak{P}$ is a path in $[\mathfrak{P}]_{s}$ satisfying that 
$$
\left(
\mathfrak{Q}, t
\right)
\leq_{\mathbf{Pth}_{\boldsymbol{\mathcal{A}}}}
\left(
\mathfrak{P}, s
\right).
$$

This completes the Case~(1).

If~(2), that is, if $\mathfrak{P}$ is a path of length strictly greater than one  containing at least one echelon and its first echelon occurs at position $i\in\bb{\mathfrak{P}}$, we consider the case (2.1) $i=0$ and the subcases (2.1.1) $\mathfrak{Q}=\mathfrak{P}^{0,0}$ and (2.1.2) $\mathfrak{Q}=\mathfrak{P}^{1,\bb{\mathfrak{P}}-1}$ and the case (2.2) $i>0$ and the subcases (2.2.1) $\mathfrak{Q}=\mathfrak{P}^{0,i-1}$ and (2.2.2) $\mathfrak{Q}=\mathfrak{P}^{i,\bb{\mathfrak{P}}-1}$.
Let us note that in each case, we have that $t=s$.

If~(2.1.1), i.e., if $\mathfrak{P}$ is a path of length strictly greater than one containing its first echelon on its first step and $\mathfrak{Q}=\mathfrak{P}^{0,0}$ then, since $\mathfrak{Q}$ is an echelon, we have, by Proposition~\ref{PCHEch}, that if $\mathfrak{Q}'\in[\mathfrak{Q}]_{s}$, then $\mathfrak{Q}=\mathfrak{Q}'$. Then,  we have that $\mathfrak{P}$ is a path in $[\mathfrak{P}]_{s}$ satisfying that 
$$
(\mathfrak{Q}',s)
\prec_{\mathbf{Pth}_{\boldsymbol{\mathcal{A}}}}
(\mathfrak{P},s).
$$

The case $i=0$ and $\mathfrak{Q}=\mathfrak{P}^{0,0}$ follows.

If~(2.1.2), i.e., if $\mathfrak{P}$ is a path of length strictly greater than one containing its first echelon on its first step and $\mathfrak{Q}=\mathfrak{P}^{1,\bb{\mathfrak{P}}-1}$ then, following Definition~\ref{DCH}, the value of the Curry Howard mapping at $\mathfrak{P}$ is given by
$$
\mathrm{CH}^{(1)}_{s}\left(
\mathfrak{P}\right)=
\mathrm{CH}^{(1)}_{s}\left(
\mathfrak{Q}
\right)
\circ_{s}^{0\mathbf{T}_{\Sigma^{\boldsymbol{\mathcal{A}}}}(X)}
\mathrm{CH}^{(1)}_{s}\left(
\mathfrak{P}^{0,0}
\right).
$$

Thus, if $\mathfrak{Q}'$ is a path in $[\mathfrak{Q}]_{s}$, then in virtue of Lemma~\ref{LCH}, 
$
\mathrm{sc}_{s}^{(0,1)}(\mathfrak{Q}')
=
\mathrm{sc}_{s}^{(0,1)}(\mathfrak{Q}).
$
Moreover,  the length of $\mathfrak{Q}'$ equals the length of $\mathfrak{Q}$.
Therefore, since $\mathfrak{P}$ decomposes as the $0$-composition $\mathfrak{P}=\mathfrak{Q}\circ^{0\mathbf{Pth}_{\boldsymbol{\mathcal{A}}}}_{s}\mathfrak{P}^{0,0}$, we have that $\mathrm{sc}^{(0,1)}_{s}(\mathfrak{Q})=\mathrm{tg}^{(0,1)}_{s}(\mathfrak{P}^{0,0})$. All in all, we conclude that $\mathrm{sc}^{(0,1)}_{s}(\mathfrak{Q}')=\mathrm{tg}^{(0,1)}_{s}(\mathfrak{P}^{0,0})$. Therefore the path $\mathfrak{P}'$ defined by
$$
\mathfrak{P}'=
\mathfrak{Q}'
\circ^{0\mathbf{Pth}_{\boldsymbol{\mathcal{A}}}}_{s}
\mathfrak{P}^{0,0}
$$
is a path of length strictly greater than one containing an echelon on its first step. Taking into account Definition~\ref{DOrd}, we have that $(\mathfrak{Q}',s)\prec_{\mathbf{Pth}_{\boldsymbol{\mathcal{A}}}} (\mathfrak{P}',s)$. Moreover the value of the Curry-Howard mapping at $\mathfrak{P}'$ is given by
\begin{align*}
\mathrm{CH}^{(1)}_{s}\left(\mathfrak{P}'\right)&=
\mathrm{CH}^{(1)}_{s}\left(
\mathfrak{Q}'
\right)
\circ_{s}^{0\mathbf{T}_{\Sigma^{\boldsymbol{\mathcal{A}}}}(X)}
\mathrm{CH}^{(1)}_{s}\left(
\mathfrak{P}^{0,0}
\right)
\\&=
\mathrm{CH}^{(1)}_{s}\left(
\mathfrak{Q}
\right)
\circ_{s}^{0\mathbf{T}_{\Sigma^{\boldsymbol{\mathcal{A}}}}(X)}
\mathrm{CH}^{(1)}_{s}
\left(
\mathfrak{P}^{0,0}
\right)
\\&=
\mathrm{CH}^{(1)}_{s}\left(
\mathfrak{P}
\right).
\end{align*}

The case $i=0$ and $\mathfrak{Q}=\mathfrak{P}^{1,\bb{\mathfrak{P}}-1}$ follows.

If~(2.2.1), i.e., if $\mathfrak{P}$ is a path of length strictly greater than one containing its first echelon on a step $i\in\bb{\mathfrak{P}}$ different from the initial one and $\mathfrak{Q}=\mathfrak{P}^{0,i-1}$ then, following Definition~\ref{DCH}, the value of the Curry-Howard mapping at $\mathfrak{P}$ is given by
$$
\mathrm{CH}^{(1)}_{s}\left(
\mathfrak{P}
\right)
=
\mathrm{CH}^{(1)}_{s}\left(
\mathfrak{P}^{i,\bb{\mathfrak{P}}-1}
\right)
\circ^{0\mathbf{T}_{\Sigma^{\boldsymbol{\mathcal{A}}}}(X)}_{s}
\mathrm{CH}^{(1)}_{s}\left(
\mathfrak{Q}
\right).
$$

Since $\mathfrak{Q}$ is an echelonless path, then in virtue of Lemma~\ref{LCHNEch} we have that
$$
\mathrm{CH}^{(1)}_{s}\left(
\mathfrak{Q}
\right)
\in \mathrm{T}_{\Sigma^{\boldsymbol{\mathcal{A}}}}(X)^{\mathsf{E}}_{s}
$$
Thus, if $\mathfrak{Q}'$ is any path in $[\mathfrak{Q}]_{s}$ then, in virtue of Lemma~\ref{LCHNEch} we conclude that $\mathfrak{Q}'$ is an echelonless path. Moreover, in virtue of Lemma~\ref{LCH}, we have that $\mathrm{tg}^{(1,0)}_{s}(\mathfrak{Q}')=\mathrm{tg}^{(1,0)}_{s}(\mathfrak{Q})$ and the length of $\mathfrak{Q}'$ equals the length of $\mathfrak{Q}$. Therefore, since $\mathfrak{P}$ decomposes as the $0$-composition $\mathfrak{P}=\mathfrak{P}^{i,\bb{\mathfrak{P}}-1}\circ^{0\mathbf{Pth}_{\boldsymbol{\mathcal{A}}}}_{s}\mathfrak{Q}$, we have that $\mathrm{sc}^{(0,1)}_{s}(\mathfrak{P}^{i,\bb{\mathfrak{P}}-1})=\mathrm{tg}^{(0,1)}_{s}(\mathfrak{Q})$. All in all, we conclude that $\mathrm{sc}^{(0,1)}_{s}(\mathfrak{P}^{i,\bb{\mathfrak{P}}-1})=\mathrm{tg}^{(0,1)}_{s}(\mathfrak{Q}')$. Therefore the path $\mathfrak{P}'$ defined by
$$
\mathfrak{P}'
=
\mathfrak{P}^{i,\bb{\mathfrak{P}}-1}
\circ_{s}^{0\mathbf{Pth}_{\boldsymbol{\mathcal{A}}}}
\mathfrak{Q}'
$$
is a path of length strictly greater than one containing its first echelon at position $i$. Taking into account Definition~\ref{DOrd}, we have that $(\mathfrak{Q}',s)\prec_{\mathbf{Pth}_{\boldsymbol{\mathcal{A}}}} (\mathfrak{P},s)$. Moreover, the value of the Curry-Howard mapping at $\mathfrak{P}'$ is given by
\begin{align*}
\mathrm{CH}^{(1)}_{s}\left(\mathfrak{P}'\right)&=
\mathrm{CH}^{(1)}_{s}\left(
\mathfrak{P}^{i,\bb{\mathfrak{P}}-1}
\right)
\circ_{s}^{0\mathbf{T}_{\Sigma^{\boldsymbol{\mathcal{A}}}}(X)}
\mathrm{CH}^{(1)}_{s}\left(
\mathfrak{Q}'
\right)
\\&=
\mathrm{CH}^{(1)}_{s}\left(
\mathfrak{P}^{i,\bb{\mathfrak{P}}-1}
\right)
\circ_{s}^{0\mathbf{T}_{\Sigma^{\boldsymbol{\mathcal{A}}}}(X)}
\mathrm{CH}^{(1)}_{s}\left(
\mathfrak{Q}
\right)
\\&=
\mathrm{CH}^{(1)}_{s}\left(
\mathfrak{P}
\right).
\end{align*}

The case $i>0$ and $\mathfrak{Q}=\mathfrak{P}^{0,i-1}$ follows.

If~(2.2.2), i.e., if $\mathfrak{P}$ is a path of length strictly greater than one containing its first echelon on a step $i\in\bb{\mathfrak{P}}$ different from the initial one and $\mathfrak{Q}=\mathfrak{P}^{i,\bb{\mathfrak{P}}-1}$ then, following Definition~\ref{DCH}, the value of the Curry-Howard mapping at $\mathfrak{P}$ is given by
$$
\mathrm{CH}^{(1)}_{s}\left(
\mathfrak{P}
\right)
=
\mathrm{CH}^{(1)}_{s}\left(
\mathfrak{Q}
\right)
\circ_{s}^{0\mathbf{T}_{\Sigma^{\boldsymbol{\mathcal{A}}}}(X)}
\mathrm{CH}^{(1)}_{s}\left(
\mathfrak{P}^{0,i-1}
\right).
$$

Since $\mathfrak{Q}$ is a path of length at least one containing an echelon on its first step, then in virtue of Lemmas~\ref{LCHEch} and~\ref{LCHEchInt} we have that
$$
\mathrm{CH}^{(1)}_{s}\left(\mathfrak{Q}\right)\in
\eta^{(1,\mathcal{A})}\left[
\mathcal{A}
\right]_{s}
\cup
 \eta^{(1,\mathcal{A})}\left[
 \mathcal{A}
 \right]^{\mathrm{int}}_{s}.
 $$
Thus, if $\mathfrak{Q}'$ is a path of length at least one containing an echelon on its first step, then in virtue of Lemmas~\ref{LCHEch} or~\ref{LCHEchInt} we conclude that $\mathfrak{Q}'$  is a path of length at least one containing an echelon on its first step. Moreover, in virtue of Lemma~\ref{LCH}, we have that $\mathrm{sc}^{(0,1)}_{s}(\mathfrak{Q}')=\mathrm{sc}^{(0,1)}_{s}(\mathfrak{Q})$ and the length of $\mathfrak{Q}'$ equals the length of $\mathfrak{Q}$. Therefore, since $\mathfrak{P}$ decomposes as the $0$-composition $\mathfrak{P}=\mathfrak{Q}\circ^{0\mathbf{Pth}_{\boldsymbol{\mathcal{A}}}}_{s}\mathfrak{P}^{0,i-1}$, we have that $\mathrm{sc}^{(0,1)}_{s}(\mathfrak{Q})=\mathrm{tg}^{(0,1)}_{s}(\mathfrak{P}^{0,i-1})$. Therefore, the path $\mathfrak{P}'$ defined by
$$
\mathfrak{P}'=
\mathfrak{Q}'
\circ^{0\mathbf{Pth}_{\boldsymbol{\mathcal{A}}}}_{s}
\mathfrak{P}^{0,i-1}
$$
is a path of length strictly greater than one containing its first echelon at position $i$. Taking into account Definition~\ref{DOrd}, we have that $(\mathfrak{Q}',s)\prec_{\mathbf{Pth}_{\boldsymbol{\mathcal{A}}}}(\mathfrak{P}',s)$. Moreover, the value of the Curry-Howard mapping at $\mathfrak{P}'$ is given by
\begin{align*}
\mathrm{CH}^{(1)}_{s}\left(
\mathfrak{P}'
\right)
&=
\mathrm{CH}^{(1)}_{s}\left(
\mathfrak{Q}'
\right)
\circ^{0\mathbf{T}_{\Sigma^{\boldsymbol{\mathcal{A}}}}(X)}_{s}
\mathrm{CH}^{(1)}_{s}\left(
\mathfrak{P}^{0,i-1}
\right)
\\&=
\mathrm{CH}^{(1)}_{s}\left(
\mathfrak{Q}
\right)
\circ^{0\mathbf{T}_{\Sigma^{\boldsymbol{\mathcal{A}}}}(X)}_{s}
\mathrm{CH}^{(1)}_{s}\left(
\mathfrak{P}^{0,i-1}
\right)
\\&=
\mathrm{CH}^{(1)}_{s}\left(
\mathfrak{P}
\right).
\end{align*}

The case $i>0$ and $\mathfrak{Q}=\mathfrak{P}^{i,\bb{\mathfrak{P}}-1}$ follows.

This completes the Case~(2).

Finally, consider the case~(3), in which $\mathfrak{P}$ is an echelonless path. Then, following Lemma~\ref{LPthHeadCt} there exists a unique word $\mathbf{s}\in S^{\star}-\{\lambda\}$ and a unique operation symbol $\sigma\in\Sigma_{\mathbf{s},s}$ associated to $\mathfrak{P}$. Let $(\mathfrak{P}_{j})_{j\in\bb{\mathbf{s}}}$ be the family of paths we can extract from $\mathfrak{P}$ in virtue of Lemma~\ref{LPthExtract}. Assume that, for an index $j\in\bb{\mathbf{s}}$, it is the case that $\mathfrak{Q}=\mathfrak{P}_{j}$. 

Then, following Definition~\ref{DCH}, the value of the Curry-Howard mapping at $\mathfrak{P}$ is given by
$$
\mathrm{CH}^{(1)}_{s}\left(
\mathfrak{P}
\right)=
\sigma^{\mathbf{T}_{\Sigma^{\boldsymbol{\mathcal{A}}}}(X)}
\left(\left(
\mathrm{CH}^{(1)}_{s_{j}}\left(
\mathfrak{P}_{j}
\right)\right)_{j\in\bb{\mathbf{s}}}
\right).
$$

Thus, if $\mathfrak{Q}'$ is a path in the path class $[\mathfrak{Q}]_{s}$ then, in virtue of Lemma~\ref{LCH}, we have that $\mathfrak{Q}'$ and $\mathfrak{Q}$ have the same length, $(0,1)$-source and $(0,1)$-target. 

Consider the path $\mathfrak{P}'$ in $\mathrm{Pth}_{\boldsymbol{\mathcal{A}},s}$ given by 
$$
\mathfrak{P}'=\sigma^{\mathbf{Pth}_{\boldsymbol{\mathcal{A}}}}
\left(
\mathfrak{P}_{0},\cdots,\mathfrak{Q}',\cdots, \mathfrak{P}_{\bb{\mathbf{s}}-1}
\right),
$$
that is, $\mathfrak{P}'$ is defined as the operation $\sigma^{\mathbf{Pth}_{\boldsymbol{\mathcal{A}}}}$ applied to the family of paths extracted from $\mathfrak{P}$, where we have replaced the path at position $j$ by $\mathfrak{Q}'$. 

Note that in virtue of Corollary~\ref{CPthWB} the path $\mathfrak{P}'$ is an echelonless path. Moreover, in virtue of Proposition~\ref{PRecov}, the path extraction algorithm from Lemma~\ref{LPthExtract} applied to $\mathfrak{P}'$ retrieves the family of paths $(
\mathfrak{P}_{0},\cdots,\mathfrak{Q}',\cdots, \mathfrak{P}_{\bb{\mathbf{s}}-1})$.

Therefore, following Definition~\ref{DOrd}, we conclude that 
$
(\mathfrak{Q}',s_{j})
\prec_{\mathbf{Pth}_{\boldsymbol{\mathcal{A}}}}
(\mathfrak{P}',s).
$ Moreover, the value of the Curry-Howard mapping at $\mathfrak{P}'$ is given by 
\allowdisplaybreaks
\begin{align*}
\mathrm{CH}^{(1)}_{s}\left(
\mathfrak{P}'
\right)&=
\sigma^{\mathbf{T}_{\Sigma^{\boldsymbol{\mathcal{A}}}}(X)}
\left(
\mathrm{CH}^{(1)}_{s_{0}}\left(\mathfrak{P}_{0}
\right),
\cdots,
\mathrm{CH}^{(1)}_{s_{j}}
\left(
\mathfrak{Q}'
\right)
\cdots,
\mathrm{CH}^{(1)}_{s_{0}}\left(
\mathfrak{P}_{0}
\right),
\right)
\\&=
\sigma^{\mathbf{T}_{\Sigma^{\boldsymbol{\mathcal{A}}}}(X)}\left(\left(
\mathrm{CH}^{(1)}_{s_{j}}\left(\mathfrak{P}_{j}
\right)\right)_{j\in\bb{\mathbf{s}}}
\right).
\\&=
\mathrm{CH}^{(1)}_{s}
\left(
\mathfrak{P}
\right).
\end{align*}

This completes the case (3) of $\mathfrak{P}$ being an echelonless path.

This concludes the base case.

\textsf{Inductive step of the induction.}

Assume the statement holds for sequences of length up to $m$, i.e., for every pair of sorts $t,s\in S$, if $\mathfrak{Q}',\mathfrak{Q}$ are paths in $\mathrm{Pth}_{\boldsymbol{\mathcal{A}},t}$ and $\mathfrak{P}$ is a path in $\mathrm{Pth}_{\boldsymbol{\mathcal{A}},s}$ such that  there exists a word $\mathbf{w}\in S^{\star}$ of length $\bb{\mathbf{w}}=m+1$ and a family of paths $(\mathfrak{R}_{k})_{k\in\bb{\mathbf{w}}}$ in $\mathrm{Pth}_{\boldsymbol{\mathcal{A}},\mathbf{w}}$  such that $w_{0}=t$, $\mathfrak{R}_{0}=\mathfrak{Q}$, $w_{m}=s$, $\mathfrak{R}_{m}=\mathfrak{P}$ and 
for every $k\in m$, $(\mathfrak{R}_{k}, w_{k})\prec_{\mathbf{Pth}_{\boldsymbol{\mathcal{A}}}} (\mathfrak{R}_{k+1}, w_{k+1})$ and  $\mathfrak{Q}'$ is a path in $[\mathfrak{Q}]_{t}$, then there exists a path $\mathfrak{P}'$ in $[\mathfrak{P}]_{s}$ such that 
$$
\left(
\mathfrak{Q}',t
\right)\leq_{\mathbf{Pth}_{\boldsymbol{\mathcal{A}}}}
\left(
\mathfrak{P}',s
\right).
$$

Now we prove it for sequences of length $m+1$. 

Let $t,s$ be sorts in $S$, if $\mathfrak{Q}',\mathfrak{Q}$ are paths in $\mathrm{Pth}_{\boldsymbol{\mathcal{A}},t}$ and $\mathfrak{P}$ is a path in $\mathrm{Pth}_{\boldsymbol{\mathcal{A}},s}$ such that  there exists a word $\mathbf{w}\in S^{\star}$ of length $\bb{\mathbf{w}}=m+2$ and a family of paths $(\mathfrak{R}_{k})_{k\in\bb{\mathbf{w}}}$ in $\mathrm{Pth}_{\boldsymbol{\mathcal{A}},\mathbf{w}}$  such that $w_{0}=t$, $\mathfrak{R}_{0}=\mathfrak{Q}$, $w_{m+1}=s$, $\mathfrak{R}_{m+1}=\mathfrak{P}$ and 
for every $k\in m+1$, $(\mathfrak{R}_{k}, w_{k})\prec_{\mathbf{Pth}_{\boldsymbol{\mathcal{A}}}} (\mathfrak{R}_{k+1}, w_{k+1})$ and let $\mathfrak{Q}'$ be a path in $[\mathfrak{Q}]_{t}$.

Consider the word $\mathbf{w}^{0,m}$ of length $\bb{\mathbf{w}^{0,m}}=m+1$ and the family of paths $(\mathfrak{R}_{k})_{k\in\bb{\mathbf{w}^{0,m}}}$ in $\mathrm{Pth}_{\boldsymbol{\mathcal{A}},\mathbf{w}^{0,m}}$. This is a sequence of length $m$ instantiating that 
$
(\mathfrak{Q},t)
\leq_{\mathbf{Pth}_{\boldsymbol{\mathcal{A}}}}
(\mathfrak{R}_{m}, w_{m}).
$ Since $\mathfrak{Q}'$ is a path in $[\mathfrak{Q}]_{t}$, we have, by the inductive hypothesis, that there exists a path $\mathfrak{R}'_{m}$ in $[\mathfrak{R}_{m}]_{w_{m}}$ satisfying that 
$
(\mathfrak{Q}',t)
\leq_{\mathbf{Pth}_{\boldsymbol{\mathcal{A}}}}
(\mathfrak{R}'_{m}, w_{m}).
$

Now, consider the sequence of paths $(\mathfrak{R}_{m},\mathfrak{P})$. It is a one-step sequence of paths in $\mathrm{Pth}_{\boldsymbol{\mathcal{A}},\mathbf{w}^{m,m+1}}$ satisfying that $(\mathfrak{R}_{m},w_{m})\prec_{\mathbf{Pth}_{\boldsymbol{\mathcal{A}}}} (\mathfrak{P}, s)$. Since $\mathfrak{R}'_{m}$ in $[\mathfrak{R}_{m}]_{w_{m}}$, by the base case,  we can find a path $\mathfrak{P}'$ in $[\mathfrak{P}]_{s}$ satisfying that 
$
(\mathfrak{R}'_{m},w_{m})
\leq_{\mathbf{Pth}_{\boldsymbol{\mathcal{A}}}}
(\mathfrak{P}', s).
$

Hence, $\mathfrak{P}'$ is a path in $[\mathfrak{P}]_{s}$ satisfying that 
$$
\left(
\mathfrak{Q}',t
\right)
\leq_{\mathbf{Pth}_{\boldsymbol{\mathcal{A}}}}
\left(
\mathfrak{P}', s
\right).
$$

This concludes the proof.
\end{proof}

\begin{corollary}\label{CCHOrdI} 
Let $t,s$ be sorts in $S$, $[\mathfrak{Q}]_{t}$ a path class in $[\mathrm{Pth}_{\boldsymbol{\mathcal{A}}}]_{t}$ and $[\mathfrak{P}]_{s}$ a path class in $[\mathrm{Pth}_{\boldsymbol{\mathcal{A}}}]_{s}$. Then the following statements are equivalent
\begin{enumerate}
\item[(i)] $([\mathfrak{Q}]_{t},t)
\leq_{[\mathbf{Pth}_{\boldsymbol{\mathcal{A}}}]}
([\mathfrak{P}]_{s},s);
$
\item[(ii)] There exists $
\mathfrak{P}'\in [\mathfrak{P}]_{s}$ such that $
(\mathfrak{Q},t)
\leq_{\mathbf{Pth}_{\boldsymbol{\mathcal{A}}}}
(\mathfrak{P}',s)
.
$
\end{enumerate}
\end{corollary}
\begin{proof}
Assume that 
$([\mathfrak{Q}]_{t},t)
\leq_{[\mathbf{Pth}_{\boldsymbol{\mathcal{A}}}]}
([\mathfrak{P}]_{s},s)
$, then there exists $\mathfrak{Q}'\in[\mathfrak{Q}]_{t}$ and $\mathfrak{P}'\in [\mathfrak{P}]_{s}$ such that $(\mathfrak{Q}',t) \leq_{\mathbf{Pth}_{\boldsymbol{\mathcal{A}}}} (\mathfrak{P}',s)$. Since $\mathfrak{Q}\in[\mathfrak{Q}']_{t}$, by Lemma~\ref{LCHOrdI},  we can find $\mathfrak{P}''\in[\mathfrak{P}']_{s}$ such that $(\mathfrak{Q},t) \leq_{\mathbf{Pth}_{\boldsymbol{\mathcal{A}}}} (\mathfrak{P}'',s)$. Note that $\mathfrak{P}''\in[\mathfrak{P}]_{s}$.

The other implication follows by definition of the relation $\leq_{[\mathbf{Pth}_{\boldsymbol{\mathcal{A}}}]}$.

This finishes the proof.
\end{proof}

In the following lemma we show that, if we are given two paths in $\mathrm{Ker}(\mathrm{CH}^{(1)})$ and one of these paths is greater than another path with respect to the order $\leq_{\mathbf{Pth}_{\boldsymbol{\mathcal{A}}}}$, then we can complete this inequality on the other side with an equivalent path.  

\begin{lemma}\label{LCHOrdII} 
Let $t,s$ be sorts in $S$, $,\mathfrak{Q}$ a path in $\mathrm{Pth}_{\boldsymbol{\mathcal{A}},t}$ and $\mathfrak{P}',\mathfrak{P}$ paths in $\mathrm{Pth}_{\boldsymbol{\mathcal{A}},s}$. If $(\mathfrak{Q},t)\leq_{\mathbf{Pth}_{\boldsymbol{\mathcal{A}}}}(\mathfrak{P},s)$ and $\mathfrak{P}'\in[\mathfrak{P}]_{s}$ then there exists a path $\mathfrak{Q}'\in [\mathfrak{Q}]_{t}$ satisfying that  
$$\left(
\mathfrak{Q}',t
\right)\leq_{\mathbf{Pth}_{\boldsymbol{\mathcal{A}}}}
\left(
\mathfrak{P}',s
\right).$$
\end{lemma}

\begin{proof}
Let us recall from Remark~\ref{ROrd} that $(\mathfrak{Q},t)\leq_{\mathbf{Pth}_{\boldsymbol{\mathcal{A}}}}(\mathfrak{P},s)$ if and only if $s=t$ and $\mathfrak{Q}=\mathfrak{P}$ or there exists a natural number $m\in\mathbb{N}-\{0\}$, a word $\mathbf{w}\in S^{\star}$ of length $\bb{\mathbf{w}}=m+1$, and a family of paths $(\mathfrak{R}_{k})$ in $\mathrm{Pth}_{\boldsymbol{\mathcal{A}},\mathbf{w}}$ such that $w_{0}=t$, $\mathfrak{R}_{0}=\mathfrak{Q}$, $w_{m}=s$, $\mathfrak{R}_{m}=\mathfrak{P}$ and for every $k\in m$, $(\mathfrak{R}_{k}, w_{k})\prec_{\mathbf{Pth}_{\boldsymbol{\mathcal{A}}}} (\mathfrak{R}_{k+1}, w_{k+1})$.

The lemma holds trivially in case $s=t$ and $\mathfrak{Q}=\mathfrak{P}$. Therefore, it remains to prove the other case. 

We will prove the lemma by induction on $m\in\mathbb{N}-\{0\}$.

\textsf{Base step of the induction.}

For $m=1$ we have that $(\mathfrak{Q},t)\prec_{\mathbf{Pth}_{\boldsymbol{\mathcal{A}}}}(\mathfrak{P},s)$, hence according to Definition~\ref{DOrd} we are in one of the following situations
\begin{enumerate}
\item $\mathfrak{P}$ and $\mathfrak{Q}$ are $(1,0)$-identity paths of the form
\begin{align*}
\mathfrak{P}&=
\mathrm{ip}^{(1,0)\sharp}_{s}(P),
&
\mathfrak{Q}&=
\mathrm{ip}^{(1,0)\sharp}_{t}(Q),
\end{align*}
for some terms $P\in\mathrm{T}_{\Sigma}(X)_{s}$ and $Q\in\mathrm{T}_{\Sigma}(X)_{t}$ and the following inequality holds
$$
\left(
Q, t
\right)
<_{\mathbf{T}_{\Sigma}(X)}
\left(
P,s
\right),
$$
where $\leq_{\mathbf{T}_{\Sigma}(X)}$ is the Artinian partial order on $\coprod\mathrm{T}_{\Sigma}(X)$ introduced in Remark~\ref{RTermOrd}.

\item $\mathfrak{P}$ is a path of length strictly greater than one  containing at least one echelon, and if its first echelon occurs at position $i\in\bb{\mathfrak{P}}$, then
\subitem if $i=0$, then $\mathfrak{Q}$ is equal to $\mathfrak{P}^{0,0}$ or $\mathfrak{P}^{1,\bb{\mathfrak{P}}-1}$,
\subitem if $i>0$, then $\mathfrak{Q}$ is equal to $\mathfrak{P}^{0,i-1}$ or $\mathfrak{P}^{i,\bb{\mathfrak{P}}-1}$;
\item $\mathfrak{P}$ is an echelonless path and $\mathfrak{Q}$ is one of the paths we can extract from $\mathfrak{P}$ in virtue of Lemma~\ref{LPthExtract}.
\end{enumerate}

If~(1), we note that both $\mathfrak{P}$ and $\mathfrak{Q}$ are $(1,0)$-identity paths. Then, if $\mathfrak{P}'$ is a path in $[\mathfrak{P}]_{s}$ then we conclude, in virtue of Corollary~\ref{CCHUZId}, that $\mathfrak{P}'=\mathfrak{P}$. This case follows easily, since $\mathfrak{Q}$ is a path in $[\mathfrak{Q}]_{t}$ satisfying that
$$
\left(
\mathfrak{Q},t
\right)
\leq_{\mathbf{Pth}_{\boldsymbol{\mathcal{A}}}}
\left(
\mathfrak{P},s
\right).
$$

This completes the case~(1).

If~(2), that is, if $\mathfrak{P}$ is a path of length strictly greater than one  containing at least one echelon and its first echelon occurs at position $i\in\bb{\mathfrak{P}}$, we consider the case~(2.1) $i=0$ and the case~(2.2) $i>0$.

If~(2.1), i.e., if $\mathfrak{P}$ is a path of length strictly greater than one containing its first echelon on its first step then, following Definition~\ref{DCH}, the value of the Curry-Howard mapping at $\mathfrak{P}$ is given by
$$
\mathrm{CH}^{(1)}_{s}\left(
\mathfrak{P}
\right)
=
\mathrm{CH}^{(1)}_{s}\left(
\mathfrak{P}^{1,\bb{\mathfrak{P}}-1}
\right)
\circ^{0\mathbf{T}_{\Sigma^{\boldsymbol{\mathcal{A}}}}(X)}_{s}
\mathrm{CH}^{(1)}_{s}\left(
\mathfrak{P}^{0,0}
\right).
$$

Since $\mathfrak{P}$ is a path of length strictly greater than one containing an echelon on its first step then, in virtue of Lemma~\ref{LCHEchInt}, we have that $\mathrm{CH}^{(1)}_{s}(\mathfrak{P})\in\eta^{(1,\mathcal{A})}[\mathcal{A}]^{\mathrm{int}}_{s}$. Thus, if $\mathfrak{P}'$ is any path in $[\mathfrak{P}]_{s}$ then, in virtue of Lemma~\ref{LCHEchInt}, we have that $\mathfrak{P}'$ is a path of length strictly greater than one containing an echelon on its first step.

Thus, the value of the Curry-Howard mapping at $\mathfrak{P}'$ is given by
$$
\mathrm{CH}^{(1)}_{s}\left(
\mathfrak{P}'
\right)
=
\mathrm{CH}^{(1)}_{s}\left(
\mathfrak{P}'^{1,\bb{\mathfrak{P}'}-1}
\right)
\circ^{0\mathbf{T}_{\Sigma^{\boldsymbol{\mathcal{A}}}}(X)}_{s}
\mathrm{CH}^{(1)}_{s}\left(
\mathfrak{P}'^{0,0}
\right).
$$

Now, we consider the different possibilities for $\mathfrak{Q}$. Note that since $(\mathfrak{Q},t)\prec_{\mathbf{Pth}_{\boldsymbol{\mathcal{A}}}}(\mathfrak{P},s)$, then following Definition~\ref{DOrd}, either (2.1.1) $\mathfrak{Q}=\mathfrak{P}^{0,0}$ or (2.1.2) $\mathfrak{Q}=\mathfrak{P}^{1,\bb{\mathfrak{P}}-1}$. In any case, note that $t=s$.

If~(2.1.1), i.e., if $\mathfrak{Q}=\mathfrak{P}^{0,0}$ then, following Definition~\ref{DOrd}, 
$(\mathfrak{P}'^{0,0},s)\prec_{\mathbf{Pth}_{\boldsymbol{\mathcal{A}}}}(\mathfrak{P}',s)$. Moreover, since $\mathfrak{P}'$ is a path in $[\mathfrak{P}]_{s}$ we conclude that $\mathrm{CH}^{(1)}_{s}(\mathfrak{Q})=\mathrm{CH}^{(1)}_{s}(\mathfrak{P}'^{0,0})$.

On the other hand, if (2.1.2), i.e., if $\mathfrak{Q}=\mathfrak{P}^{1,\bb{\mathfrak{P}}-1}$ then, following Definition~\ref{DOrd}, 
$(\mathfrak{P}'^{1,\bb{\mathfrak{P}'}-1},s)\prec_{\mathbf{Pth}_{\boldsymbol{\mathcal{A}}}}(\mathfrak{P}',s)$. Moreover, since $\mathfrak{P}'$ is a path in $[\mathfrak{P}]_{s}$ we conclude that $\mathrm{CH}^{(1)}_{s}(\mathfrak{Q})=\mathrm{CH}^{(1)}_{s}(\mathfrak{P}'^{1,\bb{\mathfrak{P}'}-1})$.

The case $i=0$ follows.

If~(2.2), i.e., if $\mathfrak{P}$ is a path of length strictly greater than one containing its first echelon on the step $i\in\bb{\mathfrak{P}}$ and $i>0$, then following Definition~\ref{DCH}, the value of the Curry-Howard mapping at $\mathfrak{P}$ is given by
$$
\mathrm{CH}^{(1)}_{s}\left(
\mathfrak{P}
\right)
=
\mathrm{CH}^{(1)}_{s}\left(
\mathfrak{P}^{i,\bb{\mathfrak{P}}-1}
\right)
\circ^{0\mathbf{T}_{\Sigma^{\boldsymbol{\mathcal{A}}}}(X)}_{s}
\mathrm{CH}^{(1)}_{s}\left(
\mathfrak{P}^{0,i-1}
\right).
$$

Since $\mathfrak{P}$ is a path of length strictly greater than one containing an echelon on a step different from the initial one then, in virtue of Lemma~\ref{LCHEchNInt}, we have that $\mathrm{CH}^{(1)}_{s}(\mathfrak{P})\in\eta^{(1,\mathcal{A})}[\mathcal{A}]^{\neg\mathrm{int}}_{s}$. Thus, if $\mathfrak{P}'$ is any path in $[\mathfrak{P}]_{s}$ then, in virtue of Lemma~\ref{LCHEchNInt}, we have that $\mathfrak{P}'$ is a path of length strictly greater than one containing an echelon on a step different from the initial one. Let $j\in\bb{\mathfrak{P}'}$ be the first index for which $\mathfrak{P}'^{j,j}$ is an echelon.

Thus, the value of the Curry-Howard mapping at $\mathfrak{P}'$ is given by
$$
\mathrm{CH}^{(1)}_{s}\left(
\mathfrak{P}'
\right)
=
\mathrm{CH}^{(1)}_{s}\left(
\mathfrak{P}'^{j,\bb{\mathfrak{P}'}-1}
\right)
\circ^{0\mathbf{T}_{\Sigma^{\boldsymbol{\mathcal{A}}}}(X)}_{s}
\mathrm{CH}^{(1)}_{s}\left(
\mathfrak{P}'^{0,j-1}
\right).
$$

Now, we consider the different possibilities for $\mathfrak{Q}$. Note that since $(\mathfrak{Q},t)\prec_{\mathbf{Pth}_{\boldsymbol{\mathcal{A}}}}(\mathfrak{P},s)$, then following Definition~\ref{DOrd}, either (2.2.1) $\mathfrak{Q}=\mathfrak{P}^{0,i-1}$ or (2.2.2) $\mathfrak{Q}=\mathfrak{P}^{i,\bb{\mathfrak{P}}-1}$. In any case, note that $t=s$.

If~(2.2.1), i.e., if $\mathfrak{Q}=\mathfrak{P}^{0,i-1}$ then, following Definition~\ref{DOrd}, 
$(\mathfrak{P}'^{0,j-1},s)\prec_{\mathbf{Pth}_{\boldsymbol{\mathcal{A}}}}(\mathfrak{P}',s)$. Moreover, since $\mathfrak{P}'$ is a path in $[\mathfrak{P}]_{s}$ we conclude that $\mathrm{CH}^{(1)}_{s}(\mathfrak{Q})=\mathrm{CH}^{(1)}_{s}(\mathfrak{P}'^{0,j-1})$.

On the other hand, if (2.2.2), i.e., if $\mathfrak{Q}=\mathfrak{P}^{j,\bb{\mathfrak{P}}-1}$ then, following Definition~\ref{DOrd}, 
$(\mathfrak{P}'^{j,\bb{\mathfrak{P}'}-1},s)\prec_{\mathbf{Pth}_{\boldsymbol{\mathcal{A}}}}(\mathfrak{P}',s)$. Moreover, since $\mathfrak{P}'$ is a path in $[\mathfrak{P}]_{s}$ we conclude that $\mathrm{CH}^{(1)}_{s}(\mathfrak{Q})=\mathrm{CH}^{(1)}_{s}(\mathfrak{P}'^{j,\bb{\mathfrak{P}'}-1})$.

The case $i>0$ follows.

This completes the Case~(2).

Finally, consider the case(3), in which $\mathfrak{P}$ is an echelonless path. 

Then, following Lemma~\ref{LPthHeadCt} there exists a unique word $\mathbf{s}\in S^{\star}-\{\lambda\}$ and a unique operation symbol $\sigma\in\Sigma_{\mathbf{s},s}$ associated to $\mathfrak{P}$. Let $(\mathfrak{P}_{j})_{j\in\bb{\mathbf{s}}}$ be the family of paths we can extract from $\mathfrak{P}$ in virtue of Lemma~\ref{LPthExtract}. Then, following Definition~\ref{DCH}, the value of the Curry-Howard mapping at $\mathfrak{P}$ is given by
$$
\mathrm{CH}^{(1)}_{s}\left(
\mathfrak{P}
\right)
=
\sigma^{\mathbf{T}_{\Sigma^{\boldsymbol{\mathcal{A}}}}(X)}\left(\left(
\mathrm{CH}^{(1)}_{s_{j}}\left(
\mathfrak{P}_{j}
\right)\right)_{j\in\bb{\mathbf{s}}}\right).
$$

Since $\mathfrak{P}$ is an echelonless path associated to the operation symbol $\sigma$ then, following Lemma~\ref{LCHNEch}, we have that $\mathrm{CH}^{(1)}_{s}(\mathfrak{P})\in\mathcal{T}(\sigma,\mathrm{T}_{\Sigma^{\boldsymbol{\mathcal{A}}}}(X))_{1}$, which is a subset of $[\mathrm{T}_{\Sigma^{\boldsymbol{\mathcal{A}}}}(X)]^{\mathsf{E}}_{s}$. Thus, if $\mathfrak{P}'$ is any path in $[\mathfrak{P}]_{s}$ then, in virtue of Lemma~\ref{LCHNEch}, we have that $\mathfrak{P}'$ is an echelonless path associated to the operation symbol $\sigma$.

Let $(\mathfrak{P}'_{j})_{j\in\bb{\mathbf{s}}}$ be the family of paths we can extract from $\mathfrak{P}$ in virtue of Lemma~\ref{LPthExtract}. Thus, the value of the Curry-Howard mapping at $\mathfrak{P}$ is given by
$$
\mathrm{CH}^{(1)}_{s}\left(
\mathfrak{P}'
\right)
=
\sigma^{\mathbf{T}_{\Sigma^{\boldsymbol{\mathcal{A}}}}(X)}\left(\left(
\mathrm{CH}^{(1)}_{s_{j}}\left(
\mathfrak{P}'_{j}
\right)\right)_{j\in\bb{\mathbf{s}}}\right).
$$

Now, we consider the different possibilities for $\mathfrak{Q}$. Note that, since $(\mathfrak{Q},t)\prec_{\mathbf{Pth}_{\boldsymbol{\mathcal{A}}}} (\mathfrak{P},s)$ then, following Definition~\ref{DOrd}, there exists an index $j\in\bb{\mathbf{s}}$ for which $\mathfrak{Q}=\mathfrak{P}_{j}$. In this case $t=s_{j}$.

Then, following Definition~\ref{DOrd}, $(\mathfrak{P}'_{j},s_{j})\prec_{\mathbf{Pth}_{\boldsymbol{\mathcal{A}}}} (\mathfrak{P}',s)$. Moreover, since $\mathfrak{P}'$ is a path in $[\mathfrak{P}]_{s}$ we conclude that $\mathrm{CH}^{(1)}_{s_{j}}(\mathfrak{Q})=\mathrm{CH}^{(1)}_{s_{j}}(\mathfrak{P}'_{j})$.

This completes the Case~(3).

This concludes the base step.

\textsf{Inductive step of the induction.}

Assume the statement holds for sequences of length up to $m$, i.e., for every pair of sorts $t,s\in S$, if $\mathfrak{Q}$ is a path in $\mathrm{Pth}_{\boldsymbol{\mathcal{A}},t}$ and $\mathfrak{P}',\mathfrak{P}$ are paths in $\mathrm{Pth}_{\boldsymbol{\mathcal{A}},s}$ such that  there exists a word $\mathbf{w}\in S^{\star}$ of length $\bb{\mathbf{w}}=m+1$ and a family of paths $(\mathfrak{R}_{k})_{k\in\bb{\mathbf{w}}}$ in $\mathrm{Pth}_{\boldsymbol{\mathcal{A}},\mathbf{w}}$  such that $w_{0}=t$, $\mathfrak{R}_{0}=\mathfrak{Q}$, $w_{m}=s$, $\mathfrak{R}_{m}=\mathfrak{P}$ and 
for every $k\in m$, $(\mathfrak{R}_{k}, w_{k})\prec_{\mathbf{Pth}_{\boldsymbol{\mathcal{A}}}} (\mathfrak{R}_{k+1}, w_{k+1})$ and  $\mathfrak{P}'$ is a path in $[\mathfrak{P}]_{s}$, then there exists a path $\mathfrak{Q}'$ in $[\mathfrak{P}]_{t}$ such that 
$$
\left(
\mathfrak{Q}',t
\right)\leq_{\mathbf{Pth}_{\boldsymbol{\mathcal{A}}}}
\left(
\mathfrak{P}',s
\right).
$$

Now we prove it for sequences of length $m+1$. 

Let $t,s$ be sorts in $S$, if $\mathfrak{Q}$ is a path in $\mathrm{Pth}_{\boldsymbol{\mathcal{A}},t}$ and $\mathfrak{P}',\mathfrak{P}$ are paths in $\mathrm{Pth}_{\boldsymbol{\mathcal{A}},s}$ such that  there exists a word $\mathbf{w}\in S^{\star}$ of length $\bb{\mathbf{w}}=m+2$ and a family of paths $(\mathfrak{R}_{k})_{k\in\bb{\mathbf{w}}}$ in $\mathrm{Pth}_{\boldsymbol{\mathcal{A}},\mathbf{w}}$  such that $w_{0}=t$, $\mathfrak{R}_{0}=\mathfrak{Q}$, $w_{m+1}=s$, $\mathfrak{R}_{m+1}=\mathfrak{P}$ and 
for every $k\in m+1$, $(\mathfrak{R}_{k}, w_{k})\prec_{\mathbf{Pth}_{\boldsymbol{\mathcal{A}}}} (\mathfrak{R}_{k+1}, w_{k+1})$ and let  $\mathfrak{P}'$ be a path in $[\mathfrak{P}]_{s}$.

Consider the word $\mathbf{w}^{1,m+1}$ of length $\bb{\mathbf{w}^{1,m+1}}=m+1$ and the family of paths $(\mathfrak{R}_{k})_{k\in\bb{\mathbf{w}^{1,m+1}}}$ in $\mathrm{Pth}_{\boldsymbol{\mathcal{A}},\mathbf{w}^{1,m+1}}$. This is a sequence of length $m$ instantiating that 
$
(\mathfrak{R}_{1},w_{1})
\leq_{\mathbf{Pth}_{\boldsymbol{\mathcal{A}}}}
(\mathfrak{P}, s).
$ Since $\mathfrak{P}'$ is a path in $[\mathfrak{P}]_{s}$, by the inductive hypothesis, there exists a path $\mathfrak{R}'_{1}$ in $[\mathfrak{R}_{1}]_{w_{1}}$ satisfying that 
$
(\mathfrak{R}'_{1},t)
\leq_{\mathbf{Pth}_{\boldsymbol{\mathcal{A}}}}
(\mathfrak{P}',s).
$

Now, consider the sequence of paths $(\mathfrak{Q},\mathfrak{R}_{1})$. It is a one-step sequence of paths in $\mathrm{Pth}_{\boldsymbol{\mathcal{A}},\mathbf{w}^{0,1}}$ satisfying that $(\mathfrak{Q},t)\prec_{\mathbf{Pth}_{\boldsymbol{\mathcal{A}}}} (\mathfrak{R}_{1}, w_{1})$. Since $\mathfrak{R}'_{1}$ in $[\mathfrak{R}_{1}]_{w_{1}}$, we have, by the base case, that we can find a path $\mathfrak{Q}'$ in $[\mathfrak{Q}]_{s}$ satisfying that 
$
(\mathfrak{Q}',t)
\leq_{\mathbf{Pth}_{\boldsymbol{\mathcal{A}}}}
(\mathfrak{R}'_{1}, w_{1}).
$

Hence, $\mathfrak{Q}'$ is a path in $[\mathfrak{Q}]_{s}$ satisfying that 
$$
\left(
\mathfrak{Q}',t
\right)
\leq_{\mathbf{Pth}_{\boldsymbol{\mathcal{A}}}}
\left(
\mathfrak{P}', s
\right).
$$

This finishes the proof.
\end{proof}

\begin{corollary}\label{CCHOrdII} 
Let $t,s$ be sorts in $S$, $[\mathfrak{Q}]_{t}$ a path class in $[\mathrm{Pth}_{\boldsymbol{\mathcal{A}}}]_{t}$ and $[\mathfrak{P}]_{s}$ a path class in $[\mathrm{Pth}_{\boldsymbol{\mathcal{A}}}]_{s}$. Then the following statements are equivalent
\begin{enumerate}
\item[(i)] $([\mathfrak{Q}]_{t},t)
\leq_{[\mathbf{Pth}_{\boldsymbol{\mathcal{A}}}]}
([\mathfrak{P}]_{s},s);
$
\item[(ii)] There exists $
\mathfrak{Q}'\in [\mathfrak{Q}]_{t}$ such that $
(\mathfrak{Q}',t)
\leq_{\mathbf{Pth}_{\boldsymbol{\mathcal{A}}}}
(\mathfrak{P},s)
.
$
\end{enumerate}
\end{corollary}
\begin{proof}
Assume that $([\mathfrak{Q}]_{t},t)
\leq_{[\mathbf{Pth}_{\boldsymbol{\mathcal{A}}}]}
([\mathfrak{P}]_{s},s)
$ then there exists $\mathfrak{Q}'\in[\mathfrak{Q}]_{t}$ and $\mathfrak{P}'\in [\mathfrak{P}]_{s}$ such that $(\mathfrak{Q}',t) \leq_{\mathbf{Pth}_{\boldsymbol{\mathcal{A}}}} (\mathfrak{P}',s)$. Since $\mathfrak{P}\in[\mathfrak{P}']_{s}$, by Lemma~\ref{LCHOrdII}, we can find $\mathfrak{Q}''\in[\mathfrak{Q}']_{t}$ such that $(\mathfrak{Q}'',t) \leq_{\mathbf{Pth}_{\boldsymbol{\mathcal{A}}}} (\mathfrak{P},s)$. Note that $\mathfrak{Q}''\in[\mathfrak{Q}]_{t}$.

The other implication follows by definition of the relation $\leq_{[\mathbf{Pth}_{\boldsymbol{\mathcal{A}}}]}$.

This finishes the proof.
\end{proof}

In the following lemma we show that if we are given two paths in $\mathrm{Ker}(\mathrm{CH}^{(1)})$, then these paths can only be related with respect to the order $\leq_{\mathbf{Pth}_{\boldsymbol{\mathcal{A}}}}$ if they are equal.

\begin{lemma}\label{LCHOrdIII} 
Let $s$ be a sort in $S$ and $\mathfrak{P}',\mathfrak{P}$ paths in $\mathrm{Pth}_{\boldsymbol{\mathcal{A}},s}$. If $(\mathfrak{P}', \mathfrak{P})\in\mathrm{Ker}(\mathrm{CH}^{(1)})_{s}$ and $(\mathfrak{P}',s)\leq_{\mathbf{Pth}_{\boldsymbol{\mathcal{A}}}}(\mathfrak{P},s)$, then $\mathfrak{P}'=\mathfrak{P}$.
\end{lemma}
\begin{proof}
Assume towards a contradiction that $\mathfrak{P}'\neq\mathfrak{P}$. Since $(\mathfrak{P}',s)\leq_{\mathbf{Pth}_{\boldsymbol{\mathcal{A}}}}(\mathfrak{P},s)$, we conclude that $(\mathfrak{P}',s)$ is strictly smaller than $(\mathfrak{P},s)$ with respect to the partial order $\leq_{\mathbf{Pth}_{\boldsymbol{\mathcal{A}}}}$. Following Proposition~\ref{PCHMono} this entails that $(\mathrm{CH}^{(1)}_{s}(\mathfrak{P}'),s)$ is strictly smaller than $(\mathrm{CH}^{(1)}_{s}(\mathfrak{P}),s)$ with respect to the partial order $\leq_{\mathbf{T}_{\Sigma^{\boldsymbol{\mathcal{A}}}}(X)}$, contradicting the fact that, by hypothesis, $(\mathfrak{P}',\mathfrak{P})\in\mathrm{Ker}(\mathrm{CH}^{(1)})_{s}$. Therefore, $\mathfrak{P}'=\mathfrak{P}$.
\end{proof}

As a corollary of the just stated lemma we have that no pair can be found with respect to the partial order $\leq_{\mathbf{Pth}_{\boldsymbol{\mathcal{A}}}}$ between two equivalent paths.

\begin{corollary}\label{CCHOrdIII} 
Let $s,t$ be sorts in $S$, $\mathfrak{Q}$ a path in $\mathrm{Pth}_{\boldsymbol{\mathcal{A}},t}$ and $\mathfrak{P}',\mathfrak{P}$ paths in $\mathrm{Pth}_{\boldsymbol{\mathcal{A}},s}$. If 
$$
\left(
\mathfrak{P}',s
\right)
\leq_{\mathbf{Pth}_{\boldsymbol{\mathcal{A}}}}
\left(
\mathfrak{Q},t
\right)
\leq_{\mathbf{Pth}_{\boldsymbol{\mathcal{A}}}}
\left(
\mathfrak{P},s
\right)
$$
and  $(\mathfrak{P}',\mathfrak{P})\in\mathrm{Ker}(\mathrm{CH}^{(1)})_{s}$, then 
$t=s$ and $\mathfrak{Q}=\mathfrak{P}$.
\end{corollary}
\begin{proof}
According to Lemma~\ref{LCHOrdIII}, $\mathfrak{P}'=\mathfrak{P}$. The statement follows from the antisymmetry of the preorder $\leq_{\mathbf{Pth}_{\boldsymbol{\mathcal{A}}}}$ proven in Proposition~\ref{POrdArt}.
\end{proof}

We are now in a position to prove that $\leq_{[\mathbf{Pth}_{\boldsymbol{\mathcal{A}}}]}$ is an Artinian order on 
$\coprod[\mathrm{Pth}_{\boldsymbol{\mathcal{A}}}]$.

\begin{restatable}{proposition}{PCHOrdArt}
\label{PCHOrdArt}  $(\coprod[\mathrm{Pth}_{\boldsymbol{\mathcal{A}}}], \leq_{[\mathbf{Pth}_{\boldsymbol{\mathcal{A}}}]})$ is an ordered set. Moreover, there is not any strictly decreasing $\omega_{0}$-chain in it, i.e., 
 $(\coprod[\mathrm{Pth}_{\boldsymbol{\mathcal{A}}}], \leq_{[\mathbf{Pth}_{\boldsymbol{\mathcal{A}}}]})$ is an Artinian  ordered set.
\end{restatable}
\begin{proof}
That $\leq_{[\mathbf{Pth}_{\boldsymbol{\mathcal{A}}}]}$ is reflexive follows from the fact that $\leq_{\mathbf{Pth}_{\boldsymbol{\mathcal{A}}}}$ is reflexive. Indeed, let $s$ be a sort in $S$ and $[\mathfrak{P}]_{s}$ a path class in $\mathrm{Pth}_{\boldsymbol{\mathcal{A}},s}$. Then, since $\mathfrak{P}\in [\mathfrak{P}]_{s}$ and $(\mathfrak{P},s)\leq_{\mathbf{Pth}_{\boldsymbol{\mathcal{A}}}} (\mathfrak{P},s)$, we have that 
$$
\left(\left[
\mathfrak{P}
\right]_{s},s\right)
\leq_{[\mathbf{Pth}_{\boldsymbol{\mathcal{A}}}]} 
\left(\left[
\mathfrak{P}\right]_{s},s\right).$$

We now prove that $\leq_{[\mathbf{Pth}_{\boldsymbol{\mathcal{A}}}]}$ is antisymmetric. Let $t,s$ be sorts in $S$ and let $[\mathfrak{Q}]_{t}$ be a path class in $[\mathrm{Pth}_{\boldsymbol{\mathcal{A}}}]_{t}$ and let $[\mathfrak{P}]_{s}$ be a path class in $[\mathrm{Pth}_{\boldsymbol{\mathcal{A}}}]_{s}$ satisfying that
\begin{itemize}
\item[(i)] $([\mathfrak{Q}]_{t},t)\leq_{[\mathbf{Pth}_{\boldsymbol{\mathcal{A}}}]}([\mathfrak{P}]_{s},s);$
\item[(ii)] $([\mathfrak{P}]_{s},s)\leq_{[\mathbf{Pth}_{\boldsymbol{\mathcal{A}}}]}([\mathfrak{Q}]_{t},t).
$
\end{itemize}

We want to prove that $t=s$ and $[\mathfrak{Q}]_{t}=[\mathfrak{P}]_{s}$.

On one hand, since the inequality (i) holds, we can find paths $\mathfrak{Q}'\in[\mathfrak{Q}]_{t}$ and $\mathfrak{P}'\in[\mathfrak{P}]_{s}$ satisfying that 
$
(\mathfrak{Q}',t)\leq_{\mathbf{Pth}_{\boldsymbol{\mathcal{A}}}}
(\mathfrak{P}',s).
$ On the other hand, since the inequality (ii) holds, we can find paths  $\mathfrak{P}''\in[\mathfrak{P}]_{s}$ and $\mathfrak{Q}''\in[\mathfrak{Q}]_{t}$ satisfying that 
$
(\mathfrak{P}'',s)\leq_{\mathbf{Pth}_{\boldsymbol{\mathcal{A}}}}
(\mathfrak{Q}'',t).
$

Since $
(\mathfrak{Q}',t)\leq_{\mathbf{Pth}_{\boldsymbol{\mathcal{A}}}}
(\mathfrak{P}',s)
$ and $\mathfrak{Q}''$ is a path in $[\mathfrak{Q}']_{t}$, in virtue of Lemma~\ref{LCHOrdI},  we can find a path $\mathfrak{P}'''\in[\mathfrak{P}']_{s}$ for which
$
(\mathfrak{Q}'',t)\leq_{\mathbf{Pth}_{\boldsymbol{\mathcal{A}}}}
(\mathfrak{P}''',s).
$

Hence, we have found the chain
$$
\left(
\mathfrak{P}'',s
\right)
\leq_{\mathbf{Pth}_{\boldsymbol{\mathcal{A}}}}
\left(
\mathfrak{Q}'',t
\right)
\leq_{\mathbf{Pth}_{\boldsymbol{\mathcal{A}}}}
\left(
\mathfrak{P}''',s
\right).
$$

Moreover, $(\mathfrak{P}'',\mathfrak{P}''')\in\mathrm{Ker}(\mathrm{CH}^{(1)})_{s}$. Hence, in virtue of Corollary~\ref{CCHOrdIII}, we conclude that $t=s$ and $\mathfrak{Q}=\mathfrak{P}$, which implies that $[\mathfrak{Q}]_{t}=[\mathfrak{P}]_{s}$.

We now prove that $\leq_{[\mathbf{Pth}_{\boldsymbol{\mathcal{A}}}]}$ is transitive. Let $u,t,s$ be sorts in $S$ and let $[\mathfrak{R}]_{u}$ be a path class $[\mathrm{Pth}_{\boldsymbol{\mathcal{A}}}]_{u}$,  $[\mathfrak{Q}]_{t}$ be a path class in $[\mathrm{Pth}_{\boldsymbol{\mathcal{A}}}]_{t}$ and let $[\mathfrak{P}]_{s}$ be a path class in $[\mathrm{Pth}_{\boldsymbol{\mathcal{A}}}]_{s}$ satisfying that
\begin{itemize}
\item[(i)] $([\mathfrak{R}]_{u},u)\leq_{[\mathbf{Pth}_{\boldsymbol{\mathcal{A}}}]}([\mathfrak{Q}]_{t},t);$
\item[(ii)] $([\mathfrak{Q}]_{t},t)\leq_{[\mathbf{Pth}_{\boldsymbol{\mathcal{A}}}]}([\mathfrak{P}]_{s},s).$
\end{itemize}

We want to prove that $([\mathfrak{R}]_{u},u)\leq_{[\mathbf{Pth}_{\boldsymbol{\mathcal{A}}}]}([\mathfrak{P}]_{s},s)$.

On one hand, since inequality (i) holds, we can find paths $\mathfrak{R}'\in[\mathfrak{R}]_{u}$ and $\mathfrak{Q}'\in[\mathfrak{Q}]_{t}$ satisfying that 
$
(\mathfrak{R}',u)\leq_{\mathbf{Pth}_{\boldsymbol{\mathcal{A}}}}
(\mathfrak{Q}',t).
$ On the other hand, since inequality (ii) holds, we can find paths  $\mathfrak{Q}''\in[\mathfrak{Q}]_{t}$ and $\mathfrak{P}'\in[\mathfrak{P}]_{s}$ satisfying that 
$
(\mathfrak{Q}'',t)\leq_{\mathbf{Pth}_{\boldsymbol{\mathcal{A}}}}
(\mathfrak{P}',s).
$

Since $
(\mathfrak{Q}'',t)\leq_{\mathbf{Pth}_{\boldsymbol{\mathcal{A}}}}
(\mathfrak{P}',s)
$ and $\mathfrak{Q}'$ is a path in $[\mathfrak{Q}'']_{t}$,  in virtue of Lemma~\ref{LCHOrdI}, we can find a path $\mathfrak{P}''\in[\mathfrak{P}']_{s}$ for which
$
(\mathfrak{Q}'',t)\leq_{\mathbf{Pth}_{\boldsymbol{\mathcal{A}}}}
(\mathfrak{P}'',s).
$

Hence, we have found the chain
$$
\left(
\mathfrak{R}',u
\right)
\leq_{\mathbf{Pth}_{\boldsymbol{\mathcal{A}}}}
\left(
\mathfrak{Q}'',t
\right)
\leq_{\mathbf{Pth}_{\boldsymbol{\mathcal{A}}}}
\left(
\mathfrak{P}'',s
\right).
$$

By the transitivity of $\leq_{\mathbf{Pth}_{\boldsymbol{\mathcal{A}}}}$, we conclude that $(\mathfrak{R}',u)
\leq_{\mathbf{Pth}_{\boldsymbol{\mathcal{A}}}}
(\mathfrak{P}'',s)$. Note that $\mathfrak{R}'\in[\mathfrak{R}]_{u}$ and $\mathfrak{P}''\in[\mathfrak{P}]_{s}$. Hence, we have that 
$$
\left(\left[
\mathfrak{R}
\right]_{u},u\right)
\leq_{[\mathbf{Pth}_{\boldsymbol{\mathcal{A}}}]}
\left(\left[
\mathfrak{P}
\right]_{s},s
\right).$$

From this it follows that $\leq_{[\mathbf{Pth}_{\boldsymbol{\mathcal{A}}}]}$ is a partial order.

We next prove that no strictly decreasing $\omega_{0}$-chain can exist in the partially ordered set  $(\coprod[\mathrm{Pth}_{\boldsymbol{\mathcal{A}}}],\leq_{[\mathbf{Pth}_{\boldsymbol{\mathcal{A}}}]})$. Let us suppose, towards a contradiction, that there exists one strictly decreasing $\omega_{0}$-chain $(([\mathfrak{R}_{k}]_{s_{k}},s_{k}))_{k\in\omega_{0}}$ in $(\coprod[\mathrm{Pth}_{\boldsymbol{\mathcal{A}}}],\leq_{[\mathbf{Pth}_{\boldsymbol{\mathcal{A}}}]})$. 

Then, for every $k\in\omega_{0}$, we have that
\begin{enumerate}
\item[(i)] $[\mathfrak{R}_{k}]_{s_{k}}$ is a path class in $[\mathrm{Pth}_{\boldsymbol{\mathcal{A}}}]_{s_{k}}$, and
\item[(ii)] $([\mathfrak{R}_{k+1}]_{s_{k+1}}, s_{k+1})<_{[\mathbf{Pth}_{\boldsymbol{\mathcal{A}}}]} ([\mathfrak{R}_{k}]_{s_{k}}, s_{k})$.
\end{enumerate}

We next prove how this strictly decreasing $\omega_{0}$-chain in the partially ordered set $(\coprod[\mathrm{Pth}_{\boldsymbol{\mathcal{A}}}],\leq_{[\mathbf{Pth}_{\boldsymbol{\mathcal{A}}}]})$ leads to a strictly decreasing $\omega_{0}$-chain in the partially ordered set $(\coprod\mathrm{Pth}_{\boldsymbol{\mathcal{A}}},\leq_{\mathbf{Pth}_{\boldsymbol{\mathcal{A}}}})$. 

Let us define $\mathfrak{R}'_{0}$ to be equal to $\mathfrak{R}_{0}$. Since $([\mathfrak{R}_{1}]_{s_{1}}, s_{1})<_{[\mathbf{Pth}_{\boldsymbol{\mathcal{A}}}]} ([\mathfrak{R}'_{0}]_{s_{0}}, s_{0})$ then, in virtue of Corollary~\ref{CCHOrdII}, we can find $\mathfrak{R}'_{1}\in[\mathfrak{R}_{1}]_{s_{1}}$ satisfying that $(\mathfrak{R}'_{1},s_{1})$ $<_{\mathbf{Pth}_{\boldsymbol{\mathcal{A}}}}$-precedes $(\mathfrak{R}'_{0},s_{0})$. 

Following this procedure,  assume that, for $n\in\mathbb{N}$ with $n\geq 2$, we have found a strictly-decreasing $n$-chain $((\mathfrak{R}'_{k},s_{k}))_{k\in n}$ in $(\coprod\mathrm{Pth}_{\boldsymbol{\mathcal{A}}},\leq_{\mathbf{Pth}_{\boldsymbol{\mathcal{A}}}})$ satisfying that, for every $k\in n$, $\mathfrak{R}'_{k}\in[\mathfrak{R}_{k}]_{s_{k}}$. Now, since $([\mathfrak{R}_{n}]_{s_{n}}, s_{n})<_{[\mathbf{Pth}_{\boldsymbol{\mathcal{A}}}]} ([\mathfrak{R}_{n-1}]_{s_{n-1}}, s_{n-1})$ and $\mathfrak{R}'_{n-1}\in [\mathfrak{R}_{n-1}]_{s_{n-1}}$ then, in virtue of Corollary~\ref{CCHOrdII},  we can find $\mathfrak{R}'_{n}\in[\mathfrak{R}_{n}]_{s_{n}}$ satisfying that $(\mathfrak{R}'_{n},s_{n})<_{\mathbf{Pth}_{\boldsymbol{\mathcal{A}}}}(\mathfrak{R}'_{n-1},s_{n-1})$.  Note that  $((\mathfrak{R}'_{k},s_{k}))_{k\in n+1}$ is  a strictly-decreasing $n+1$-chain in $(\coprod\mathrm{Pth}_{\boldsymbol{\mathcal{A}}},\leq_{\mathbf{Pth}_{\boldsymbol{\mathcal{A}}}})$ satisfying that, for every $k\in n+1$, $\mathfrak{R}'_{k}\in[\mathfrak{R}_{k}]_{s_{k}}$.

This will ultimately leads to a strictly-decreasing $\omega_{0}$-chain in $(\coprod\mathrm{Pth}_{\boldsymbol{\mathcal{A}}},\leq_{\mathbf{Pth}_{\boldsymbol{\mathcal{A}}}})$ contradicting the fact that, by Proposition~\ref{POrdArt}, $(\coprod\mathrm{Pth}_{\boldsymbol{\mathcal{A}}},\leq_{\mathbf{Pth}_{\boldsymbol{\mathcal{A}}}})$ is an Artinian partially ordered set.

It follows that $(\coprod[\mathrm{Pth}_{\boldsymbol{\mathcal{A}}}],\leq_{[\mathbf{Pth}_{\boldsymbol{\mathcal{A}}}]})$ is an Artinian partially ordered set.

This concludes the proof.
\end{proof}

We next show that several mappings from and to $\coprod[\mathrm{Pth}_{\boldsymbol{\mathcal{A}}}]$ are order-preserving, order-reflecting or order-embeddings.

\begin{restatable}{proposition}{PCHOrd}
\label{PCHOrd} 
The mapping $\coprod\mathrm{pr}^{\mathrm{Ker}(\mathrm{CH}^{(1)})}\colon\coprod\mathrm{Pth}_{\boldsymbol{\mathcal{A}}}\mor\coprod[\mathrm{Pth}_{\boldsymbol{\mathcal{A}}}]$ that, for every $s\in S$, maps a pair $(\mathfrak{P},s)$ in $\coprod\mathrm{Pth}_{\boldsymbol{\mathcal{A}}}$ to the pair $([\mathfrak{P}]_{s},s)$ in $\coprod[\mathrm{Pth}_{\boldsymbol{\mathcal{A}}}]$ is order-preserving and inversely order-preserving.
\end{restatable}

Taking into account the definitions above, now the coproduct of the $(1,0)$-identity path mapping becomes an order embedding from the Artinian partial order $(\coprod\mathrm{T}_{\Sigma}(X), \leq_{\mathbf{T}_{\Sigma}(X)})$, introduced on Remark~\ref{RTermOrd}, to the Artinian partial order $(\coprod[\mathrm{Pth}_{\boldsymbol{\mathcal{A}}}], \leq_{[\mathbf{Pth}_{\boldsymbol{\mathcal{A}}}]})$, introduced on Definition~\ref{DCHOrd}.

\begin{restatable}{proposition}{PCHOrdIp}
\label{PCHOrdIp} 
The mapping $\coprod\mathrm{ip}^{([1],0)\sharp}$ is an order-embedding
$$
\textstyle
\coprod\mathrm{ip}^{([1],0)\sharp}\colon
\left(
\coprod\mathrm{T}_{\Sigma}(X)
\leq_{\mathbf{T}_{\Sigma}(X)}
\right)
\mor 
\left(\coprod[\mathrm{Pth}_{\boldsymbol{\mathcal{A}}}], \leq_{[\mathbf{Pth}_{\boldsymbol{\mathcal{A}}}]}
\right).
$$
\end{restatable}
\begin{proof}
Let $s,t$ be sorts in $S$ and $(Q,t), (P,s)$ pairs in $\coprod\mathrm{T}_{\Sigma}(X)$. We need to prove that the following statements are equivalent
\begin{enumerate}
\item[(i)] $(Q,t)\leq_{\mathbf{T}_{\Sigma}(X)} (P,s)$;
\item[(ii)] $(\mathrm{ip}^{([1],0)\sharp}_{t}(Q),t)
\leq_{[\mathbf{Pth}_{\boldsymbol{\mathcal{A}}}]}
(\mathrm{ip}^{([1],0)\sharp}_{s}(P),s).
$ 
\end{enumerate}

Note that the following chain of equivalences holds
\allowdisplaybreaks
\begin{align*}
(Q,t)\leq_{\mathbf{T}_{\Sigma}(X)} (P,s)
&\Longleftrightarrow
\left(
\mathrm{ip}^{(0,1)\sharp}_{t}
\left(Q
\right), t
\right)
\leq_{\mathbf{Pth}_{\boldsymbol{\mathcal{A}}}}
\left(
\mathrm{ip}^{(0,1)\sharp}_{s}
\left(
P\right), s
\right)
\tag{1}
\\&\Longleftrightarrow
\left(\left[
\mathrm{ip}^{(0,1)\sharp}_{t}
\left(
Q
\right)\right]_{t}, t
\right)
\leq_{[\mathbf{Pth}_{\boldsymbol{\mathcal{A}}}]}
\left(
\left[
\mathrm{ip}^{(0,1)\sharp}_{s}
\left(
P\right)\right]_{s}, s\right)
\tag{2}
\\&\Longleftrightarrow
\left(\mathrm{ip}^{(0,[1])\sharp}_{t}\left(
Q
\right), t\right)
\leq_{[\mathbf{Pth}_{\boldsymbol{\mathcal{A}}}]}
\left(
\mathrm{ip}^{(0,[1])\sharp}_{s}
\left(
P\right), s\right).
\tag{3}
\end{align*}

The first equivalence follows from the fact that, by Proposition~\ref{POrdEmb}, the mapping $\coprod\mathrm{ip}^{(1,0)\sharp}$ is an order embedding; the second equivalence follows from left to right from the fact that, by Proposition~\ref{PCHOrd}, $\coprod\mathrm{pr}^{\mathrm{Ker}(\mathrm{CH}^{(1)})}$ is order-preserving. On the other hand the implication from right to left follows from the fact that, by Proposition~\ref{PCHOrd}, $\coprod\mathrm{pr}^{\mathrm{Ker}(\mathrm{CH}^{(1)})}$ is inversely order-preserving and from the fact that, by Corollary~\ref{CCHUZId}, the class under the kernel of the Curry-Howard mapping of a $(1,0)$-identity path only contains this exact $(1,0)$-identity path; finally, the last equivalence follows from the fact that, according to Definition~\ref{DCHUZ}, $\mathrm{ip}^{([1],0)\sharp}=\mathrm{pr}^{\mathrm{Ker}(\mathrm{CH}^{(1)})}\circ\mathrm{ip}^{(1,0)\sharp}$.

This concludes the proof.
\end{proof}

\begin{restatable}{proposition}{PCHOrdMono}
\label{PCHOrdMono} 
The mapping $\coprod\mathrm{CH}^{(1)\mathrm{m}}$ from
$\coprod[\mathrm{Pth}_{\boldsymbol{\mathcal{A}}}]$ to  $\coprod\mathrm{T}_{\Sigma^{\boldsymbol{\mathcal{A}}}(X)}$ determines an order-preserving mapping 
$$
\textstyle
\coprod\mathrm{CH}^{(1)\mathrm{m}}\colon
\left(\coprod[\mathrm{Pth}_{\boldsymbol{\mathcal{A}}}], \leq_{[\mathbf{Pth}_{\boldsymbol{\mathcal{A}}}]}
\right)
\mor
\left(\coprod\mathrm{T}_{\Sigma^{\boldsymbol{\mathcal{A}}}(X)}, \leq_{\mathbf{T}_{\Sigma^{\boldsymbol{\mathcal{A}}}}(X)}
\right),
$$
i.e., given pairs $([\mathfrak{Q}]_{t},t),([\mathfrak{P}]_{s},s)$ in $\coprod[\mathrm{Pth}_{\boldsymbol{\mathcal{A}}}]$ satisfying that
$$\left(\left[\mathfrak{Q}\right]_{t},t\right)\leq_{[\mathbf{Pth}_{\boldsymbol{\mathcal{A}}}]} \left(\left[\mathfrak{P}\right]_{s},s\right),$$
then 
$$\left(
\mathrm{CH}^{(1)\mathrm{m}}_{t}
\left(\left[
\mathfrak{Q}\right]_{t}\right),t\right)\leq_{\mathbf{T}_{\Sigma^{\boldsymbol{\mathcal{A}}}}(X)}\left(
\mathrm{CH}^{(1)\mathrm{m}}_{s}
\left(\left[
\mathfrak{P}
\right]_{s}
\right),s
\right),$$
that is, $\mathrm{CH}^{(1)}_{t}(\mathfrak{Q})$ is a subterm of type $t$ of the term $\mathrm{CH}^{(1)}_{s}(\mathfrak{P})$.
\end{restatable}
\begin{proof}
It follows from Proposition~\ref{PCHMono}
\end{proof}
\chapter{
\texorpdfstring
{The free completion of the first-order identity path mapping}
{The free completion}
}\label{S1G}

In this chapter we introduce the free completion of the partial $\Sigma^{\boldsymbol{\mathcal{A}}}$-algebra 
$\mathbf{Pth}_{\boldsymbol{\mathcal{A}}}$, denoted by $\mathbf{F}_{\Sigma^{\boldsymbol{\mathcal{A}}}}(\mathbf{Pth}_{\boldsymbol{\mathcal{A}}})$ (which is a particular case of Proposition~\ref{PFreeComp}, where we introduced the free completion of a partial $\Sigma$-algebra). This will allow us to consider the free completion of the $S$-sorted mapping $\mathrm{ip}^{(1,X)}$ from $X$ to $\mathrm{Pth}_{\boldsymbol{\mathcal{A}}}$, denoted by $\mathrm{ip}^{(1,X)@}$, which, by construction, is a $\Sigma^{\boldsymbol{\mathcal{A}}}$-homomorphism from $\mathbf{T}_{\Sigma^{\boldsymbol{\mathcal{A}}}}(X)$ to $\mathbf{F}_{\Sigma^{\boldsymbol{\mathcal{A}}}}(\mathbf{Pth}_{\boldsymbol{\mathcal{A}}})$. This mapping will allow us to show that, for every sort $s\in S$ and every path $\mathfrak{P}\in \mathrm{Pth}_{\boldsymbol{\mathcal{A}},s}$, the object $\mathrm{ip}^{(1,X)@}_{s}(\mathrm{CH}^{(1)}_{s}(\mathfrak{P}))$ which, a priori, is an object in $\mathbf{F}_{\Sigma^{\boldsymbol{\mathcal{A}}}}(\mathbf{Pth}_{\boldsymbol{\mathcal{A}}})_{s}$ is indeed a path in $\mathrm{Pth}_{\boldsymbol{\mathcal{A}},s}$. Moreover, we show that $\mathrm{ip}^{(1,X)@}_{s}(\mathrm{CH}^{(1)}_{s}(\mathfrak{P}))$ is, in fact, a path in $[\mathfrak{P}]_{s}\subseteq \mathrm{Pth}_{\boldsymbol{\mathcal{A}},s}$. This has several interesting consequences. We have that $\mathrm{ip}^{(1,X)@}_{s}(\mathrm{CH}^{(1)}_{s}(\mathfrak{P}))$ is equal to $\mathfrak{P}$, for  a $(1,0)$-identity path 
$\mathfrak{P}$. We prove that the composition $\mathrm{ip}^{(1,X)@}\circ \mathrm{CH}^{(1)}$ is a $\Sigma$-homomorphism that also preserves rewrite rules, the $0$-source and the $0$-target operation symbols. Furthermore, for a one-step path $\mathfrak{P}$ in $\mathrm{Pth}_{\boldsymbol{\mathcal{A}},s}$, we also obtain that $\mathrm{ip}^{(1,X)@}_{s}(\mathrm{CH}^{(1)}_{s}(\mathfrak{P}))$ is equal to 
$\mathfrak{P}$. We conclude this chapter by showing that the mappings $\mathrm{ip}^{(1,X)@}$ and  $\mathrm{ip}^{(1,X)@}\circ \mathrm{CH}^{(1)}$ are order-preserving.


We begin by applying the free completion construction to the partial $\Sigma^{\boldsymbol{\mathcal{A}}}$-algebra $\mathbf{Pth}_{\boldsymbol{\mathcal{A}}}$ and presenting the basic results necessary to understand its  behavior. The following definition and the statements contained in it follow directly from Proposition~\ref{PFreeComp}.

\begin{restatable}{definition}{DFP}
\label{DFP} 
\index{free completion!first-order!$\mathbf{F}_{\Sigma^{\boldsymbol{\mathcal{A}}}}(\mathbf{Pth}_{\boldsymbol{\mathcal{A}}})$}
For the partial $\Sigma^{\boldsymbol{\mathcal{A}}}$-algebra $\mathbf{Pth}_{\boldsymbol{\mathcal{A}}}$ we let 
$\mathbf{F}_{\Sigma^{\boldsymbol{\mathcal{A}}}}(\mathbf{Pth}_{\boldsymbol{\mathcal{A}}})$ stand for the free 
$\Sigma^{\boldsymbol{\mathcal{A}}}$-completion of $\mathbf{Pth}_{\boldsymbol{\mathcal{A}}}$. Moreover, 
\index{inclusion!first-order!$\eta^{(1,\mathbf{Pth}_{\boldsymbol{\mathcal{A}}})}$}
we denote by $\eta^{(1,\mathbf{Pth}_{\boldsymbol{\mathcal{A}}})}$ the $\Sigma^{\boldsymbol{\mathcal{A}}}$-homomorphism
$$
\eta^{(1,\mathbf{Pth}_{\boldsymbol{\mathcal{A}}})}\colon
\mathbf{Pth}_{\boldsymbol{\mathcal{A}}}
\mor
\mathbf{F}_{\Sigma^{\boldsymbol{\mathcal{A}}}}\left(
\mathbf{Pth}_{\boldsymbol{\mathcal{A}}}
\right)
$$
given by the insertion of generators. let us note that $\eta^{(1,\mathbf{Pth}_{\boldsymbol{\mathcal{A}}})}$ is an embedding. As usual, we identify $\mathrm{Pth}_{\boldsymbol{\mathcal{A}}}$ with its image under $\eta^{(1,\mathbf{Pth}_{\boldsymbol{\mathcal{A}}})}$. In this way, $\mathbf{Pth}_{\boldsymbol{\mathcal{A}}}$ becomes a weak subalgebra of $\mathbf{F}_{\Sigma^{\boldsymbol{\mathcal{A}}}}(\mathbf{Pth}_{\boldsymbol{\mathcal{A}}})$.
\end{restatable}

\begin{remark}\label{RFP} 
Let us recall that, by Proposition~\ref{PFreeComp}, for every $\Sigma^{\boldsymbol{\mathcal{A}}}$-algebra 
$\mathbf{B}$ and every $\Sigma^{\boldsymbol{\mathcal{A}}}$-homomorphism $f$ from $\mathbf{Pth}_{\boldsymbol{\mathcal{A}}}$ to $\mathbf{B}$, there exists a unique $\Sigma^{\boldsymbol{\mathcal{A}}}$-homomorphism $f^{\mathrm{fc}}$ from $\mathbf{F}_{\Sigma^{\boldsymbol{\mathcal{A}}}}(\mathbf{Pth}_{\boldsymbol{\mathcal{A}}})$ to $\mathbf{B}$ such that $f=f^{\mathrm{fc}}\circ\eta^{(1,\mathbf{Pth}_{\boldsymbol{\mathcal{A}}})}$.
\end{remark}

We next apply the same construction to the discrete $\Sigma^{\boldsymbol{\mathcal{A}}}$-algebra on $X$, i.e., the partial $\Sigma^{\boldsymbol{\mathcal{A}}}$-algebra whose underlying $S$-sorted set is $X$ and whose operations are nowhere defined.

\begin{definition}\label{DFX} We let $\mathbf{D}_{\Sigma^{\boldsymbol{\mathcal{A}}}}(X)$ stand for the partial $\Sigma^{\boldsymbol{\mathcal{A}}}$-algebra whose underlying $S$-sorted set is $X$ and which is such that, for each $(\mathbf{s},s)\in S^{\star}\times S$ and for each operation symbol in $\Sigma^{\boldsymbol{\mathcal{A}}}_{\mathbf{s},s}$, the domain of definition of the corresponding operation is empty. We call it the \emph{discrete} $\Sigma^{\boldsymbol{\mathcal{A}}}$-algebra \emph{on} $X$. For $\mathbf{D}_{\Sigma^{\boldsymbol{\mathcal{A}}}}(X)$ it happens that every $S$-sorted mapping $f$ from $X$ to a partial $\Sigma^{\boldsymbol{\mathcal{A}}}$-algebra $\mathbf{B}$ is a $\Sigma^{\boldsymbol{\mathcal{A}}}$-homomorphism, 
$$
f\colon\mathbf{D}_{\Sigma^{\boldsymbol{\mathcal{A}}}}(X)
\mor
\mathbf{B}.
$$
Moreover, by Remark~\ref{RFreeDisc}, $\mathbf{F}_{\Sigma^{\boldsymbol{\mathcal{A}}}}(
\mathbf{D}_{\Sigma^{\boldsymbol{\mathcal{A}}}}(X))$, the free completion of $\mathbf{D}_{\Sigma^{\boldsymbol{\mathcal{A}}}}(X)$, and $\mathbf{T}_{\Sigma^{\boldsymbol{\mathcal{A}}}}(X)$, the free $\Sigma^{\boldsymbol{\mathcal{A}}}$-algebra on $X$, are such that $\mathbf{F}_{\Sigma^{\boldsymbol{\mathcal{A}}}}(
\mathbf{D}_{\Sigma^{\boldsymbol{\mathcal{A}}}}(X))= \mathbf{T}_{\Sigma^{\boldsymbol{\mathcal{A}}}}(X)$. For this reason, the standard insertion of generators $\eta^{(1,X)}$ from $X$ to $\mathrm{T}_{\Sigma^{\boldsymbol{\mathcal{A}}}}(X)$ becomes a $\Sigma^{\boldsymbol{\mathcal{A}}}$-homomorphism
$$
\eta^{(1,X)}\colon
\mathbf{D}_{\Sigma^{\boldsymbol{\mathcal{A}}}}(X)
\mor
\mathbf{T}_{\Sigma^{\boldsymbol{\mathcal{A}}}}(X).
$$
Let us note that $\eta^{(1,X)}$ is an embedding. As usual, we identify $\mathbf{D}_{\Sigma^{\boldsymbol{\mathcal{A}}}}(X)$ with its image under $\eta^{(1,X)}$. In this way, $\mathbf{D}_{\Sigma^{\boldsymbol{\mathcal{A}}}}(X)$ becomes a weak subalgebra of $\mathbf{T}_{\Sigma^{\boldsymbol{\mathcal{A}}}}(X)$.
\end{definition}

\begin{remark}\label{RFX} 
Let us recall that, by Proposition~\ref{PFreeComp}, for every $\Sigma^{\boldsymbol{\mathcal{A}}}$-algebra 
$\mathbf{B}$ and every $\Sigma^{\boldsymbol{\mathcal{A}}}$-homomorphism $f$ from $\mathbf{D}_{\Sigma^{\boldsymbol{\mathcal{A}}}}(X)$ to $\mathbf{B}$, there exists a unique $\Sigma^{\boldsymbol{\mathcal{A}}}$-homomorphism $f^{\mathrm{fc}}$ from $\mathbf{F}_{\Sigma^{\boldsymbol{\mathcal{A}}}}(\mathbf{D}_{\Sigma^{\boldsymbol{\mathcal{A}}}}(X))$ to $\mathbf{B}$ such that $f=f^{\mathrm{fc}}\circ\eta^{(1,X)}$(let us note that this is the universal property of $\mathbf{T}_{\Sigma^{\boldsymbol{\mathcal{A}}}}(X)$). 
\end{remark}

The following definition and the statements contained in it follow directly from Propositions~\ref{PFreeComp} and~\ref{PSch}.

\begin{restatable}{definition}{DIp}
\label{DIp} 
The $S$-sorted mapping $\mathrm{ip}^{(1,X)}$ from $X$ to $\mathbf{Pth}_{\boldsymbol{\mathcal{A}}}$ that, for every  sort $s\in S$, sends a variable $x\in X_{s}$ to $\mathrm{ip}^{(1,X)}_{s}(x)=\mathrm{ip}^{(1,0)\sharp}_{s}(\eta^{(0,X)}_{s}(x))$, the $(1,0)$-identity  path on $x$, is a $\Sigma^{\boldsymbol{\mathcal{A}}}$-homomorphism from $\mathbf{D}_{\Sigma^{\boldsymbol{\mathcal{A}}}}(X)$ to $\mathbf{Pth}_{\boldsymbol{\mathcal{A}}}$,
$$
\mathrm{ip}^{(1,X)}
\colon
\mathbf{D}_{\Sigma^{\boldsymbol{\mathcal{A}}}}(X)
\mor
\mathbf{Pth}_{\boldsymbol{\mathcal{A}}}.
$$
\index{identity!first-order!$\mathrm{ip}^{(1,X)\mathrm{fc}}$}
On the other hand, by Corollary~\ref{CFreeAdj}, there exists a unique $\Sigma^{\boldsymbol{\mathcal{A}}}$-homomorphism, denoted $\mathrm{ip}^{(1,X)@}$, from $\mathbf{T}_{\Sigma^{\boldsymbol{\mathcal{A}}}}(X)$ to $\mathbf{F}_{\Sigma^{\boldsymbol{\mathcal{A}}}}(\mathbf{Pth}_{\boldsymbol{\mathcal{A}}})$,
$$
\mathrm{ip}^{(1,X)@}
\colon
\mathbf{T}_{\Sigma^{\boldsymbol{\mathcal{A}}}}(X)
\mor
\mathbf{F}_{\Sigma^{\boldsymbol{\mathcal{A}}}}(\mathbf{Pth}_{\boldsymbol{\mathcal{A}}})
$$ 
such that 
$
\mathrm{ip}^{(1,X)@}\circ
\eta^{(1,X)}
=
\eta^{(1,\mathbf{Pth}_{\boldsymbol{\mathcal{A}}})}
\circ
\mathrm{ip}^{(1,X)}.
$
We will refer to $\mathrm{ip}^{(1,X)@}$ as the \emph{free completion of the $S$-sorted mapping $\mathrm{ip}^{(1,X)}$}. 
In this particular case we have that 
$$
\mathrm{ip}^{(1,X)@}=\left(
\eta^{(1,\mathbf{Pth}_{\boldsymbol{\mathcal{A}}})}\circ \mathrm{ip}^{(1,X)}
\right)^{\mathrm{fc}}
$$
and that
$$
\mathbf{Sch}\left(\mathrm{ip}^{(1,X)}
\right)=\left(\left(\eta^{(1,\mathbf{Pth}_{\boldsymbol{\mathcal{A}}})}\circ \mathrm{ip}^{(1,X)}\right)^{\mathrm{fc}}
\right)^{-1}
\left[\eta^{(1,\mathbf{Pth}_{\boldsymbol{\mathcal{A}}})}
\left[\mathbf{Pth}_{\boldsymbol{\mathcal{A}}}
\right]\right].
$$
where, by Proposition~\ref{PSch}, $\mathbf{Sch}(\mathrm{ip}^{(1,X)})$ is the $X$-generated relative subalgebra of the free $\Sigma^{\boldsymbol{\mathcal{A}}}$-completion of $\mathbf{D}_{\Sigma^{\boldsymbol{\mathcal{A}}}}(X)$, which, as we know, is $\mathbf{T}_{\Sigma^{\boldsymbol{\mathcal{A}}}}(X)$. Finally, we let $\mathrm{ip}^{(1,X)\mathrm{Sch}}$ stand for the unique closed $\Sigma^{\boldsymbol{\mathcal{A}}}$-homomorphism from $\mathbf{Sch}(\mathrm{ip}^{(1,X)})$ to $\mathbf{Pth}_{\boldsymbol{\mathcal{A}}}$ making the square of the diagram
in Figure~\ref{FIp} commute.
\end{restatable}

\begin{figure}
\begin{tikzpicture}
[ACliment/.style={-{To [angle'=45, length=5.75pt, width=4pt, round]}
}]
\node[] (A) 	at 	(-3,1.5) 	[] 	{$\mathbf{D}_{\Sigma^{\boldsymbol{\mathcal{A}}}}(X)$};
\node[] (B) 	at 	(0,-2) 	[] 	{$\mathbf{Pth}_{\boldsymbol{\mathcal{A}}}$};
\node[]	(S)		at 	(0,0)	[]	{$\mathbf{Sch}(\mathrm{ip}^{(1,X)})$};
\node[] (FA)	at	(4,0)	[]	{$\mathbf{T}_{\Sigma^{\boldsymbol{\mathcal{A}}}}(X)$};
\node[] (FB)	at	(4,-2)	[]	{$\mathbf{F}_{\Sigma^{\boldsymbol{\mathcal{A}}}}(\mathbf{Pth}_{\boldsymbol{\mathcal{A}}})$};
\draw[ACliment, bend right]  (A) 	to node [below left]	{$\mathrm{ip}^{(1,X)}$} (B);
\draw[ACliment, bend left]  (A) 	to node [above right]	{$\eta^{(1,X)}$} (FA);
\draw[ACliment]  (A) 	to node [above right]
{$\mathrm{in}^{\mathbf{D}_{\Sigma^{\boldsymbol{\mathcal{A}}}}(X),\mathbf{Sch}(\mathrm{ip}^{(1,X)})}$} (S);
\draw[ACliment]  (B) 	to node [below]	{$\eta^{(1,\mathbf{Pth}_{\boldsymbol{\mathcal{A}}})}$} (FB);
\draw[ACliment]  (S) 	to node [above]	{$\mathrm{in}^{\mathbf{Sch}(\mathrm{ip}^{(1,X)})}$} (FA);
\draw[ACliment]  (S) 	to node [right]	{$\mathrm{ip}^{(1,X)\mathrm{Sch}}$} (B);
\draw[ACliment]  (FA) 	to node [right]	{$\mathrm{ip}^{(1,X)@}$} (FB);
\end{tikzpicture}
\caption{The free completion of $\mathrm{ip}^{(1,X)}$.}
\label{FIp}
\end{figure}


\section{
\texorpdfstring
{Interaction with the Curry-Howard mapping}
{Interaction}
}
After having introduced the free completion of the $\mathrm{ip}^{(1,X)}$ mapping, we will present a series of results that aim to better understand its behaviour with respect to those terms of  
$\mathrm{T}_{\Sigma^{\boldsymbol{\mathcal{A}}}}(X)$ that are in the image of $\mathrm{Pth}_{\boldsymbol{\mathcal{A}}}$ by means of $\mathrm{CH}^{(1)}$.

We start by investigating the compositions of $\mathrm{ip}^{(1,X)@}$ with several mappings to 
$\mathbf{T}_{\Sigma^{\boldsymbol{\mathcal{A}}}}(X)$, the domain of $\mathrm{ip}^{(1,X)@}$. But before doing it, we introduce the $\Sigma$-reduct of the free $\Sigma^{\boldsymbol{\mathcal{A}}}$-completion of the many-sorted partial $\Sigma^{\boldsymbol{\mathcal{A}}}$-algebra of  paths

\begin{definition}\label{DPFRed} Let $\mathrm{in}^{\Sigma,(1,0)}$ be the canonical embedding of $\Sigma$ into $\Sigma^{\boldsymbol{\mathcal{A}}}$. Then, by Proposition~\ref{PFunSig}, for the morphism $\mathbf{in}^{\Sigma,(1,0)}=(\mathrm{id}_{S},\mathrm{in}^{\Sigma,(1,0)})$ from $(S,\Sigma)$ to $(S,\Sigma^{\boldsymbol{\mathcal{A}}})$ and the free $\Sigma^{\boldsymbol{\mathcal{A}}}$-completion $\mathbf{F}_{\Sigma^{\boldsymbol{\mathcal{A}}}}(\mathbf{Pth}_{\boldsymbol{\mathcal{A}}})$, we will denote by $\mathbf{F}_{\Sigma^{\boldsymbol{\mathcal{A}}}}^{(1,0)}(\mathbf{Pth}_{\boldsymbol{\mathcal{A}}})$ the $\Sigma$-algebra $\mathbf{in}^{\Sigma,(0,1)}(\mathbf{F}_{\Sigma^{\boldsymbol{\mathcal{A}}}}(\mathbf{Pth}_{\boldsymbol{\mathcal{A}}})$. We will call $\mathbf{F}_{\Sigma^{\boldsymbol{\mathcal{A}}}}^{(0,1)}(\mathbf{Pth}_{\boldsymbol{\mathcal{A}}})$ the $\Sigma$-reduct of the free $\Sigma^{\boldsymbol{\mathcal{A}}}$-completion $\mathbf{F}_{\Sigma^{\boldsymbol{\mathcal{A}}}}(
\mathbf{Pth}_{\boldsymbol{\mathcal{A}}})$ of $\mathbf{Pth}_{\boldsymbol{\mathcal{A}}}$.
\end{definition}

\begin{figure}
\begin{tikzpicture}
[ACliment/.style={-{To [angle'=45, length=5.75pt, width=4pt, round]},
scale=.8
}]
\node[] (X) 	at 	(0,0) 	[] 	{$X$};
\node[] (PT) 	at 	(4,0) 	[] 	{$\mathrm{T}_{\Sigma}(X)$};
\node[] (P2)	at		(4,-2)	[]	{$\mathbf{Pth}_{\boldsymbol{\mathcal{A}}}$};
\node[] (FP2)	at		(8,-2)	 []
{$\mathbf{F}_{\Sigma^{\boldsymbol{\mathcal{A}}}}(\mathbf{Pth}_{\boldsymbol{\mathcal{A}}})$};
\node[] (T2)	at		(4,-4)	[]	{$\mathbf{T}_{\Sigma^{\boldsymbol{\mathcal{A}}}}(X)$};

\draw[ACliment]  (X) 	to node [above]	{$\eta^{(0,X)}$} (PT);
\draw[ACliment, bend right=10]  (X) 	to node [below left]	 {$\mathrm{ip}^{(1,X)}$} (P2);
\draw[ACliment]  (P2) 	to node [above]	{$\eta^{(1,\mathbf{Pth}_{\boldsymbol{\mathcal{A}}})}$} (FP2);
\draw[ACliment]  (PT) 	to node [left]	{$\mathrm{ip}^{(1,0)\sharp}$} (P2);
\draw[ACliment]  (P2) 	to node [left]	{$\mathrm{CH}^{(1)}$} (T2);
\draw[ACliment, bend right=10]  (T2) 	to node [below right]	 {$\mathrm{ip}^{(1,X)@}$} (FP2);
\end{tikzpicture}
\caption{$\mathrm{ip}^{(1,X)\mathrm{fc}}$ on $(1,0)$-identity paths.}
\label{FIpUZ}
\end{figure}

\begin{proposition}\label{PIpUZ} The following equations holds
\begin{itemize}
\item[(i)] $
\mathrm{ip}^{(1,X)@}\circ\eta^{(1,0)\sharp}
=
\eta^{(1,\mathbf{Pth}_{\boldsymbol{\mathcal{A}}})}
\circ
\mathrm{ip}^{(1,0)\sharp}
;
$
\item[(ii)] $\mathrm{ip}^{(1,X)@}
\circ
\mathrm{CH}^{(1)}
\circ
\mathrm{ip}^{(1,0)\sharp}
=
\eta^{(1,\mathbf{Pth}_{\boldsymbol{\mathcal{A}}})}
\circ
\mathrm{ip}^{(1,0)\sharp}.
$
\end{itemize}

\end{proposition}
\begin{proof}
The reader is advised to consult the diagram in Figure~\ref{FIpUZ} for a better insight of the many-sorted mappings under consideration.

We first prove that the following equation holds
$$
\mathrm{ip}^{(1,X)@}\circ\eta^{(1,0)\sharp}
=
\eta^{(1,\mathbf{Pth}_{\boldsymbol{\mathcal{A}}})}
\circ
\mathrm{ip}^{(1,0)\sharp}
.
$$

Let us note that the $S$-sorted mapping $\mathrm{ip}^{(1,X)@}\circ\eta^{(1,0)\sharp}$ is a $\Sigma$-homomorphism from $\mathbf{T}_{\Sigma}(X)$ to $\mathbf{F}^{(0,1)}_{\Sigma^{\boldsymbol{\mathcal{A}}}}(\mathbf{Pth}_{\boldsymbol{\mathcal{A}}})$ by what is stated in Definition~\ref{DIp} and by Proposition~\ref{PEmb}. Moreover, the following chain of equalities holds
\begin{align*}
\mathrm{ip}^{(1,X)@}
\circ
\eta^{(1,0)\sharp}
\circ
\eta^{(0,X)}
&=
\mathrm{ip}^{(1,X)@}
\circ
\eta^{(1,X)}
\tag{1}
\\&=
\eta^{(1,\mathbf{Pth}_{\boldsymbol{\mathcal{A}}})}
\circ
\mathrm{ip}^{(1,X)}.
\tag{2}
\end{align*}

The first equality follows from Proposition~\ref{PEmb}; while the second equality follows from what is stated in Definition~\ref{DIp}.

On the other hand the $S$-sorted mapping $\eta^{(1,\mathbf{Pth}_{\boldsymbol{\mathcal{A}}})}\circ\mathrm{ip}^{(1,0)\sharp}$ is a $\Sigma$-homomorphism from $\mathbf{T}_{\Sigma}(X)$ to $\mathbf{F}^{(0,1)}_{\Sigma^{\boldsymbol{\mathcal{A}}}}(\mathbf{Pth}_{\boldsymbol{\mathcal{A}}})$ by Remark~\ref{RFP} and Proposition~\ref{PIpHom}. Moreover, the following equality holds
$$
\eta^{(1,\mathbf{Pth}_{\boldsymbol{\mathcal{A}}})}
\circ
\mathrm{ip}^{(1,0)\sharp}
\circ
\eta^{(0,X)}
=
\eta^{(1,\mathbf{Pth}_{\boldsymbol{\mathcal{A}}})}
\circ
\mathrm{ip}^{(1,X)}.
$$
This is so, by Proposition~\ref{PBasicEq}.

Hence, by the universal property of the free $\Sigma$-algebra $\mathbf{T}_{\Sigma}(X)$, we have the desired equality:
\begin{align*}
\mathrm{ip}^{(1,X)@}\circ\eta^{(1,0)\sharp}
&=
\eta^{(1,\mathbf{Pth}_{\boldsymbol{\mathcal{A}}})}
\circ
\mathrm{ip}^{(1,0)\sharp}
.
\tag{$\star$}
\end{align*}

For the second equality we have that the following chain of equalities holds
\begin{align*}
\mathrm{ip}^{(1,X)@}
\circ
\mathrm{CH}^{(1)}
\circ
\mathrm{ip}^{(1,0)\sharp}
&=
\mathrm{ip}^{(1,X)@}
\circ
\eta^{(1,0)\sharp}
\tag{1}
\\&=
\eta^{(1,\mathbf{Pth}_{\boldsymbol{\mathcal{A}}})}
\circ
\mathrm{ip}^{(1,0)\sharp}
\tag{2}
\end{align*}

The first equality follows from Proposition~\ref{PCHId}; the last equality follows from equality~($\star$) above.

This finishes the proof.
\end{proof}

\begin{proposition}
\label{PIpEch} The following equation holds
\begin{itemize}
\item[(i)] $
\mathrm{ip}^{(1,X)@}\circ\eta^{(1,\mathcal{A})}
=
\mathrm{ech}^{(1,\mathcal{A})}.
$
\end{itemize}
\end{proposition}
\begin{proof}
Let $s$ be a sort in $S$ and let $\mathfrak{p}$ be a rewrite rule in $\mathcal{A}_{s}$.

The following chain of equalities holds
\begin{align*}
\mathrm{ip}^{(1,X)@}_{s}\left(\eta^{(1,\mathcal{A})}_{s}\left(
\mathfrak{p}
\right)\right)
&=
\mathrm{ip}^{(1,X)@}_{s}\left(
\mathfrak{p}^{\mathbf{T}_{\Sigma^{\boldsymbol{\mathcal{A}}}}(X)}
\right)
\tag{1}
\\&=
\mathfrak{p}^{
\mathbf{F}_{\Sigma^{\boldsymbol{\mathcal{A}}}}(\mathbf{Pth}_{\boldsymbol{\mathcal{A}}})
}
\tag{2}
\\&=
\mathfrak{p}^{
\mathbf{Pth}_{\boldsymbol{\mathcal{A}}}
}
\tag{3}
\\&=
\mathrm{ech}^{(1,\mathcal{A})}_{s}\left(
\mathfrak{p}
\right).
\tag{4}
\end{align*}

In the just stated chain of equalities, the first equality unravels the value of the mapping $\eta^{(1,\mathcal{A})}$ at $\mathfrak{p}$, according to Definition~\ref{DEta}; the second equality follows from the fact that $\mathrm{ip}^{(1,X)@}$ is a $\Sigma^{\boldsymbol{\mathcal{A}}}$-homomorphism, by what is stated in Definition~\ref{DIp}; the third equality follows from the fact that there is an interpretation for the constant operation symbol $\mathfrak{p}$ in the many-sorted partial $\Sigma^{\boldsymbol{\mathcal{A}}}$-algebra $\mathbf{Pth}_{\boldsymbol{\mathcal{A}}}$, according to Proposition~\ref{PPthCatAlg}, thus the interpretation of the constant operation symbol $\mathfrak{p}$ in the many-sorted  $\Sigma^{\boldsymbol{\mathcal{A}}}$-algebra $\mathbf{F}_{\Sigma^{\boldsymbol{\mathcal{A}}}}(\mathbf{Pth}_{\boldsymbol{\mathcal{A}}})$ becomes that of $\mathbf{Pth}_{\boldsymbol{\mathcal{A}}}$; finally, the last equality recovers the description of the constant operation symbol $\mathfrak{p}$ in the many-sorted partial $\Sigma^{\boldsymbol{\mathcal{A}}}$-algebra $\mathbf{Pth}_{\boldsymbol{\mathcal{A}}}$, according to Proposition~\ref{PPthCatAlg}.
\end{proof}

The following proposition is fundamental for the rest of this work. It states that $\mathrm{ip}^{(1,X)@}$ acting on the value of $\mathrm{CH}^{(1)}$ at a path $\mathfrak{P}$ is always another path, not necessarily equal to the input $\mathfrak{P}$, but which belongs to the $\mathrm{Ker}(\mathrm{CH}^{(1)})$-equivalence class of $\mathfrak{P}$. This proposition takes into account the embedding of $\mathbf{Pth}_{\boldsymbol{\mathcal{A}}}$ as a weak subalgebra of $\mathbf{F}_{\Sigma^{\boldsymbol{\mathcal{A}}}}(\mathbf{Pth}_{\boldsymbol{\mathcal{A}}})$, as stated in Definition~\ref{DFP}. 

\begin{restatable}{proposition}{PIpCH}
\label{PIpCH} The mapping 
$$
\mathrm{ip}^{(1,X)@}
\circ
\mathrm{CH}^{(1)}
\colon
\mathrm{Pth}_{\boldsymbol{\mathcal{A}}}
\mor
\mathrm{F}_{\Sigma^{\boldsymbol{\mathcal{A}}}}(
\mathbf{Pth}_{\boldsymbol{\mathcal{A}}}
)
$$
sends, for every sort $s\in S$, a path $\mathfrak{P}$ in $\mathrm{Pth}_{\boldsymbol{\mathcal{A}},s}$ to a path in the class $[\mathfrak{P}]_{s}$.
\end{restatable}
\begin{proof}
We want to show that, for every sort $s\in S$ and every path $\mathfrak{P}$ in $\mathrm{Pth}_{\boldsymbol{\mathcal{A}},s}$, the element $\mathrm{ip}^{(1,X)@}_{s}(\mathrm{CH}^{(1)}_{s}(\mathfrak{P}))$ of $\mathrm{F}_{\Sigma^{\boldsymbol{\mathcal{A}}}}(\mathbf{Pth}_{\boldsymbol{\mathcal{A}}})_{s}$ is such that
\begin{itemize}
\item[(i)] $\mathrm{ip}^{(1,X)@}_{s}(\mathrm{CH}^{(1)}_{s}(\mathfrak{P}))$ is a path in $\mathrm{Pth}_{\boldsymbol{\mathcal{A}},s}$; and
\item[(ii)] $\mathrm{ip}^{(1,X)@}_{s}(\mathrm{CH}^{(1)}_{s}(\mathfrak{P}))$ belongs to $[\mathfrak{P}]_{s}$.
\end{itemize}

In case of $\mathfrak{P}$ being a $(1,0)$-identity path, these results follow from Proposition~\ref{PIpUZ}. Note that, for every sort $s\in S$, if $\mathfrak{P}$ is a $(1,0)$-identity path, then there exists a path term $P$ in $\mathrm{T}_{\Sigma}(X)_{s}$,
for which $\mathfrak{P}=\mathrm{ip}^{(1,0)\sharp}_{s}(P)$.

The following chain of equalities holds
\allowdisplaybreaks
\begin{align*}
\mathrm{ip}^{(1,X)@}_{s}\left(\mathrm{CH}^{(1)}_{s}\left(
\mathfrak{P}
\right)\right)
&=
\mathrm{ip}^{(1,X)@}_{s}\left(\mathrm{CH}^{(1)}_{s}\left(
\mathrm{ip}^{(1,0)\sharp}_{s}\left(
P
\right)
\right)\right)
\tag{1}
\\&=
\eta^{(1,\mathbf{Pth}_{\boldsymbol{\mathcal{A}}})}_{s}\left(
\mathrm{ip}^{(1,0)\sharp}_{s}\left(P
\right)\right)
\tag{2}
\\&=
\eta^{(1,\mathbf{Pth}_{\boldsymbol{\mathcal{A}}})}_{s}\left(
\mathfrak{P}
\right)
.
\tag{3}
\end{align*}

The first equality unravels the description of $\mathfrak{P}$ as a $(1,0)$-identity path; the second equality follows from Proposition~\ref{PIpUZ}; finally, the last equality recovers the description of $\mathfrak{P}$ as a $(1,0)$-identity path.

We identify $\eta^{(1,\mathbf{Pth}_{\boldsymbol{\mathcal{A}}})}_{s}(
\mathfrak{P}
)$ with $\mathfrak{P}$. Note that  
\begin{itemize}
\item[(i)] $\mathrm{ip}^{(1,X)@}_{s}(\mathrm{CH}^{(1)}_{s}(
\mathfrak{P}
))$  is a  path in $\mathrm{Pth}_{\boldsymbol{\mathcal{A}},s}$, and
\item[(ii)] $\mathrm{ip}^{(1,X)@}_{s}(\mathrm{CH}^{(1)}_{s}(
\mathfrak{P}
))$  belongs to $[\mathfrak{P}]_{s}$.
\end{itemize}

We can, therefore, assume without loss of generality that $\mathfrak{P}$ is a  path of length at least one.

We prove the general case by Artinian induction on 
$(
\coprod\mathrm{Pth}_{\boldsymbol{\mathcal{A}}},
\leq_{\mathbf{Pth}_{\boldsymbol{\mathcal{A}}}}
)
$.

\textsf{Base step of the Artinian induction}.

Let $(\mathfrak{P},s)$ be a minimal element in $(\coprod\mathrm{Pth}_{\boldsymbol{\mathcal{A}}},\leq_{\mathbf{Pth}_{\boldsymbol{\mathcal{A}}}})
$. Then by Proposition~\ref{PMinimal}, and taking into account that we are assuming that $\mathfrak{P}$ is a non-$(1,0)$-identity path, we have that $\mathfrak{P}$ is an echelon associated to a  rewrite rule $\mathfrak{p}\in\mathcal{A}_{s}$. 

According to Definition~\ref{DCH}, the value of the  Curry-Howard mapping at $\mathfrak{P}$ is given by 
$$
\mathrm{CH}^{(1)}_{s}\left(
\mathfrak{P}
\right)
=
\mathfrak{p}^{
\mathbf{T}_{\Sigma^{\boldsymbol{\mathcal{A}}}}(X)
}.
$$

Thus we have that 
\begin{align*}
\mathrm{ip}^{(1,X)@}_{s}\left(
\mathrm{CH}^{(1)}_{s}\left(
\mathfrak{P}
\right)\right)
&=
\mathrm{ip}^{(1,X)@}_{s}\left(
\mathfrak{p}^{\mathbf{T}_{\Sigma^{\boldsymbol{\mathcal{A}}}}(X)}
\right)
\tag{1}
\\&=
\mathfrak{p}^{
\mathbf{F}_{\Sigma^{\boldsymbol{\mathcal{A}}}}\left(
\mathbf{Pth}_{\boldsymbol{\mathcal{A}}}
\right)
}
\tag{2}
\\&=
\mathfrak{p}^{
\mathbf{Pth}_{\boldsymbol{\mathcal{A}}}
}.
\tag{3}
\end{align*}

In the just stated chain of equalities, the first equality unravels the Curry-Howard mapping according to Definition~\ref{DCH}; the second equality holds because, according to Definition~\ref{DIp}, 
$\mathrm{ip}^{(1,X)@}$ is a $\Sigma^{\boldsymbol{\mathcal{A}}}$-homomorphism; and the third equality holds since, by Proposition~\ref{PPthCatAlg},  we have an interpretation of the constant symbol $\mathfrak{p}$ in the many-sorted partial $\Sigma^{\boldsymbol{\mathcal{A}}}$-algebra $\mathbf{Pth}_{\boldsymbol{\mathcal{A}}}$.

Let us recall that $\mathfrak{p}^{
\mathbf{Pth}_{\boldsymbol{\mathcal{A}}}
}$ is the echelon
$\mathrm{ech}^{(1,\mathcal{A})}_{s}(\mathfrak{p})$, i.e., the original path $\mathfrak{P}$. 

Therefore, 
\begin{itemize}
\item[(i)] $\mathrm{ip}^{(1,X)@}_{s}(\mathrm{CH}^{(1)}_{s}(
\mathfrak{P}
))$  is a  path in $\mathrm{Pth}_{\boldsymbol{\mathcal{A}},s}$, and
\item[(ii)] $\mathrm{ip}^{(1,X)@}_{s}(\mathrm{CH}^{(1)}_{s}(
\mathfrak{P}
))$  belongs to $[\mathfrak{P}]_{s}$.
\end{itemize}

This concludes the base case.

\textsf{Inductive step of the Artinian induction}.

Let $(\mathfrak{P},s)$ be a non-minimal element in $(\coprod\mathrm{Pth}_{\boldsymbol{\mathcal{A}}}, \leq_{\mathbf{Pth}_{\boldsymbol{\mathcal{A}}}})$. Let us suppose that, for every sort $t\in S$ and every path $\mathfrak{Q}$ in $\mathrm{Pth}_{\boldsymbol{\mathcal{A}},t}$, if $(\mathfrak{Q},t)<_{\mathbf{Pth}_{\boldsymbol{\mathcal{A}}}} (\mathfrak{P},s)$, then the statement holds for $\mathfrak{Q}$, i.e., 
\begin{itemize}
\item[(i)] $\mathrm{ip}^{(1,X)@}_{t}(\mathrm{CH}^{(1)}_{t}(\mathfrak{Q}))$ is a  path in $\mathrm{Pth}_{\boldsymbol{\mathcal{A}},t}$; and
\item[(ii)] $\mathrm{ip}^{(1,X)@}_{t}(\mathrm{CH}^{(1)}_{t}(\mathfrak{Q}))$ belongs to $[\mathfrak{Q}]_{t}$.
\end{itemize}

Since $(\mathfrak{P},s)$ is a non-minimal element in $(\coprod\mathrm{Pth}_{\boldsymbol{\mathcal{A}}},\leq_{\mathbf{Pth}_{\boldsymbol{\mathcal{A}}}})$ and taking into account that $\mathfrak{P}$ is a non-$(1,0)$-identity path, we have, by Lemma~\ref{LOrdI}, that $\mathfrak{P}$ is either (1) a  path of length strictly greater than one containing at least one  echelon or (2) an echelonless path.

If (1), i.e., if $\mathfrak{P}$ is a  path of length strictly greater than one containing at least one  echelon. Then, we let $i\in\bb{\mathfrak{P}}$ be the first index for which the one-step subpath $\mathfrak{P}^{i,i}$ is an echelon. We distinguish two cases accordingly, either (1.1) $i=0$ or~(1.2) $i>0$.

If~(1.1), i.e., if $\mathfrak{P}$  is a path of length strictly greater than one having its first echelon on its first step then, according to Definition~\ref{DCH}, the value of the  Curry-Howard mapping at $\mathfrak{P}$ is given by
$$
\mathrm{CH}^{(1)}_{s}\left(
\mathfrak{P}
\right)
=
\mathrm{CH}^{(1)}_{s}\left(
\mathfrak{P}^{1,\bb{\mathfrak{P}}-1}
\right)
\circ_{s}^{0\mathbf{T}_{\Sigma^{\boldsymbol{\mathcal{A}}}}(X)}
\mathrm{CH}^{(1)}_{s}\left(
\mathfrak{P}^{0,0}
\right).
$$

Therefore, the  following chain of equalities holds
\begin{flushleft}
$
\mathrm{ip}^{(1,X)@}_{s}\left(
\mathrm{CH}^{(1)}_{s}\left(
\mathfrak{P}
\right)\right)$
\allowdisplaybreaks
\begin{align*}
\quad
&=
\mathrm{ip}^{(1,X)@}_{s}\left(
\mathrm{CH}^{(1)}_{s}\left(
\mathfrak{P}^{1,\bb{\mathfrak{P}}-1}
\right)
\circ_{s}^{0\mathbf{T}_{\Sigma^{\boldsymbol{\mathcal{A}}}}(X)}
\mathrm{CH}^{(1)}_{s}\left(
\mathfrak{P}^{0,0}
\right)
\right)
\tag{1}
\\&=
\mathrm{ip}^{(1,X)@}_{s}\left(
\mathrm{CH}^{(1)}_{s}\left(
\mathfrak{P}^{1,\bb{\mathfrak{P}}-1}
\right)\right)
\circ_{s}^{0
\mathbf{F}_{\Sigma^{\boldsymbol{\mathcal{A}}}}\left(
\mathbf{Pth}_{\boldsymbol{\mathcal{A}}}
\right)}
\mathrm{ip}^{(1,X)@}_{s}\left(
\mathrm{CH}^{(1)}_{s}\left(
\mathfrak{P}^{0,0}
\right)\right)
\tag{2}
\\&=
\mathrm{ip}^{(1,X)@}_{s}\left(
\mathrm{CH}^{(1)}_{s}\left(
\mathfrak{P}^{1,\bb{\mathfrak{P}}-1}
\right)\right)
\circ_{s}^{0
\mathbf{Pth}_{\boldsymbol{\mathcal{A}}}
}
\mathrm{ip}^{(1,X)@}_{s}\left(
\mathrm{CH}^{(1)}_{s}\left(
\mathfrak{P}^{0,0}
\right)\right).
\tag{3}
\end{align*}
\end{flushleft}

In the just stated chain of equalities, the first equality recovers the value of the Curry-Howard mapping at $\mathfrak{P}$; the second equality holds since, according to Definition~\ref{DIp}, $\mathrm{ip}^{(1,X)@}$ is a many-sorted $\Sigma^{\boldsymbol{\mathcal{A}}}$-homomorphism; finally, the last equality follows taking into account that $(\mathfrak{P}^{1,\bb{\mathfrak{P}}-1},s)$ and $(\mathfrak{P}^{0,0},s)$ are strictly smaller than $(\mathfrak{P},s)$ with respect to $\prec_{\mathbf{Pth}_{\boldsymbol{\mathcal{A}}}}$. Then, by the inductive hypothesis, we have that $\mathrm{ip}^{(1,X)@}_{s}(\mathrm{CH}^{(1)}_{s}(\mathfrak{P}^{1,\bb{\mathfrak{P}}-1}))$ and $\mathrm{ip}^{(1,X)@}_{s}(\mathrm{CH}^{(1)}_{s}(\mathfrak{P}^{0,0}))$ are paths in $[\mathfrak{P}^{1,\bb{\mathfrak{P}}-1}]_{s}$ and $[\mathfrak{P}^{0,0}]_{s}$ respectively. Hence the interpretation of the $0$-composition operation symbol $\circ^{0}_{s}$ in $\mathbf{F}_{\Sigma^{\boldsymbol{\mathcal{A}}}}(\mathbf{Pth}_{\boldsymbol{\mathcal{A}}})$ becomes that of $\mathbf{Pth}_{\boldsymbol{\mathcal{A}}}$. Note that the $0$-composition is defined because  paths in the same class have the same $(0,1)$-source and $(0,1)$-target according to Lemma~\ref{LCH}.

Hence $\mathrm{ip}^{(1,X)@}_{s}(\mathrm{CH}^{(1)}_{s}(\mathfrak{P}))$ is a path.

Let us also note that, since $\mathrm{ip}^{(1,X)@}_{s}(\mathrm{CH}^{(1)}_{s}(\mathfrak{P}^{0,0}))$ is a path in $[\mathfrak{P}^{0,0}]_{s}$, we have, by Lemma~\ref{LCHUEch}, that $\mathrm{ip}^{(1,X)@}_{s}(\mathrm{CH}^{(1)}_{s}(\mathfrak{P}^{0,0}))$ is an echelon. Moreover, since the path $\mathrm{ip}^{(1,X)@}_{s}(\mathrm{CH}^{(1)}_{s}(\mathfrak{P}^{1,\bb{\mathfrak{P}}-1}))$ is in the path class $[\mathfrak{P}^{1,\bb{\mathfrak{P}}-1}]_{s}$, we have, by Lemma~\ref{LCH}, that $\mathrm{ip}^{(1,X)@}_{s}(\mathrm{CH}^{(1)}_{s}(\mathfrak{P}^{1,\bb{\mathfrak{P}}-1}))$ is a  path of length equal to $\bb{\mathfrak{P}^{1,\bb{\mathfrak{P}}-1}}$. 

All in all, we can affirm that $\mathrm{ip}^{(1,X)@}_{s}(\mathrm{CH}^{(1)}_{s}(\mathfrak{P}))$ is a path of length strictly greater than one containing an echelon on its first step. Hence, according to Definition~\ref{DCH}, the value of the Curry-Howard mapping at $\mathrm{ip}^{(1,X)@}_{s}(\mathrm{CH}^{(1)}_{s}(\mathfrak{P}))$ is given by
\begin{flushleft}
$
\mathrm{CH}^{(1)}_{s}\left(
\mathrm{ip}^{(1,X)@}_{s}\left(
\mathrm{CH}^{(1)}_{s}\left(
\mathfrak{P}
\right)\right)\right)
$
\begin{align*}
\quad&=
\mathrm{CH}^{(1)}_{s}\left(
\left(
\mathrm{ip}^{(1,X)@}_{s}\left(
\mathrm{CH}^{(1)}_{s}\left(
\mathfrak{P}
\right)\right)
\right)^{1,\bb{
\mathrm{ip}^{(1,X)@}_{s}(
\mathrm{CH}^{(1)}_{s}(
\mathfrak{P}
))
}-1}
\right)
\circ^{0\mathbf{T}_{\Sigma^{\boldsymbol{\mathcal{A}}}}(X)}_{s}
\\&\qquad\qquad\qquad\qquad\qquad\qquad\qquad\qquad
\mathrm{CH}^{(1)}_{s}\left(
\left(
\mathrm{ip}^{(1,X)@}_{s}\left(
\mathrm{CH}^{(1)}_{s}\left(
\mathfrak{P}
\right)\right)
\right)^{0,0}\right)
\tag{1}
\\&=
\mathrm{CH}^{(1)}_{s}\left(
\mathrm{ip}^{(1,X)@}_{s}\left(
\mathrm{CH}^{(1)}_{s}\left(
\mathfrak{P}^{1,\bb{\mathfrak{P}}-1}
\right)\right)\right)
\circ_{s}^{0
\mathbf{T}_{\Sigma^{\boldsymbol{\mathcal{A}}}}(X)
}
\\&\qquad\qquad\qquad\qquad\qquad\qquad\qquad\qquad
\mathrm{CH}^{(1)}_{s}\left(
\mathrm{ip}^{(1,X)@}_{s}\left(
\mathrm{CH}^{(1)}_{s}\left(
\mathfrak{P}^{0,0}
\right)\right)\right)
\tag{2}
\\&=
\mathrm{CH}^{(1)}_{s}\left(
\mathfrak{P}^{1,\bb{\mathfrak{P}}-1}
\right)
\circ_{s}^{0
\mathbf{T}_{\Sigma^{\boldsymbol{\mathcal{A}}}}(X)
}
\mathrm{CH}^{(1)}_{s}\left(
\mathfrak{P}^{0,0}
\right)
\tag{3}
\\&=
\mathrm{CH}^{(1)}_{s}\left(
\mathfrak{P}
\right).
\tag{4}
\end{align*}
\end{flushleft}

The first equality follows from the fact that $\mathrm{ip}^{(1,X)@}_{s}(\mathrm{CH}^{(1)}_{s}(\mathfrak{P}))$ is a  path of length strictly greater than one containing an echelon on its first step and Definition~\ref{DCH}; the second equality follows from the fact that $\mathrm{ip}^{(1,X)@}_{s}(\mathrm{CH}^{(1)}_{s}(\mathfrak{P}^{0,0}))$ is an echelon. Hence, from the fact that $\mathrm{ip}^{(1,X)@}_{s}(\mathrm{CH}^{(1)}_{s}(\mathfrak{P}))$ can be written as the $0$-composition
$$
\mathrm{ip}^{(1,X)@}_{s}\left(
\mathrm{CH}^{(1)}_{s}\left(
\mathfrak{P}^{1,\bb{\mathfrak{P}}-1}
\right)\right)
\circ_{s}^{0\mathbf{Pth}_{\boldsymbol{\mathcal{A}}}}
\mathrm{ip}^{(1,X)@}_{s}\left(
\mathrm{CH}^{(1)}_{s}\left(
\mathfrak{P}^{0,0}
\right)\right),
$$
we infer that 
\begin{align*}
&
\mathrm{ip}^{(1,X)@}_{s}\left(
\mathrm{CH}^{(1)}_{s}\left(
\mathfrak{P}
\right)\right)^{0,0}
=
\mathrm{ip}^{(1,X)@}_{s}\left(
\mathrm{CH}^{(1)}_{s}\left(
\mathfrak{P}^{0,0}
\right)\right),  \mbox{ and} 
\\&
\mathrm{ip}^{(1,X)@}_{s}\left(
\mathrm{CH}^{(1)}_{s}\left(
\mathfrak{P}
\right)\right)^{1,\bb{\mathrm{ip}^{(1,X)@}_{s}(
\mathrm{CH}^{(1)}_{s}(
\mathfrak{P}
))}-1}
=
\mathrm{ip}^{(1,X)@}_{s}\left(
\mathrm{CH}^{(1)}_{s}\left(
\mathfrak{P}^{1,\bb{\mathfrak{P}}-1}
\right)\right);
\end{align*}
the third equality follows from the fact that the paths $\mathrm{ip}^{(1,X)@}_{s}(\mathrm{CH}^{(1)}_{s}(\mathfrak{P}^{1,\bb{\mathfrak{P}}-1}))$ and $\mathrm{ip}^{(1,X)@}_{s}(\mathrm{CH}^{(1)}_{s}(\mathfrak{P}^{0,0}))$ are paths in $[\mathfrak{P}^{1,\bb{\mathfrak{P}}-1}]_{s}$
and $[\mathfrak{P}^{0,0}]_{s}$ respectively; finally, the last equality recovers the value of the  Curry-Howard mapping at $\mathfrak{P}$.

Therefore,
\begin{itemize}
\item[(i)] $\mathrm{ip}^{(1,X)@}_{s}(\mathrm{CH}^{(1)}_{s}(\mathfrak{P}))$ is a path in $\mathrm{Pth}_{\boldsymbol{\mathcal{A}},s}$; and
\item[(ii)] $\mathrm{ip}^{(1,X)@}_{s}(\mathrm{CH}^{(1)}_{s}(\mathfrak{P}))$ belongs to $[\mathfrak{P}]_{s}$.
\end{itemize}

The case $i=0$ follows.

If~(1.2), i.e., if $\mathfrak{P}$  is a  path of length strictly greater than one having its first  echelon at position $i\in\bb{\mathfrak{P}}$, with $i>0$, then, according to Definition~\ref{DCH}, the value of the  Curry-Howard mapping at $\mathfrak{P}$ is given by
$$
\mathrm{CH}^{(1)}_{s}\left(
\mathfrak{P}
\right)
=
\mathrm{CH}^{(1)}_{s}\left(
\mathfrak{P}^{i,\bb{\mathfrak{P}}-1}
\right)
\circ_{s}^{0\mathbf{T}_{\Sigma^{\boldsymbol{\mathcal{A}}}}(X)}
\mathrm{CH}^{(1)}_{s}\left(
\mathfrak{P}^{0,i-1}
\right).
$$

This case follows from a similar argument to that presented in the case $i=0$.

If~(2), i.e., if $\mathfrak{P}$ is an echelonless path then, according to Lemma~\ref{LPthHeadCt}, there exists a unique word $\mathbf{s}\in S^{\star}-\{\lambda\}$ and a unique operation symbol $\sigma$ in $\Sigma_{\mathbf{s},s}$ associated to $\mathfrak{P}$. Let $(\mathfrak{P}_{j})_{j\in\bb{\mathbf{s}}}$ be the family of paths we can extract from $\mathfrak{P}$ in virtue of Lemma~\ref{LPthExtract}. Then, according to Definition~\ref{DCH}, we have that the value of the  Curry-Howard mapping at $\mathfrak{P}$ is given by
$$
\mathrm{CH}^{(1)}_{s}\left(
\mathfrak{P}
\right)
=
\sigma^{\mathbf{T}_{\Sigma^{\boldsymbol{\mathcal{A}}}}(X)}
\left(\left(\mathrm{CH}^{(1)}_{s_{j}}\left(
\mathfrak{P}_{j}
\right)\right)_{j\in\bb{\mathbf{s}}}
\right).
$$

Therefore, the following chain of equalities holds
\begin{align*}
\mathrm{ip}^{(1,X)@}_{s}\left(
\mathrm{CH}^{(1)}_{s}\left(
\mathfrak{P}
\right)\right)&=
\mathrm{ip}^{(1,X)@}_{s}\left(
\sigma^{\mathbf{T}_{\Sigma^{\boldsymbol{\mathcal{A}}}}(X)}
\left(\left(\mathrm{CH}^{(1)}_{s_{j}}\left(
\mathfrak{P}_{j}
\right)\right)_{j\in\bb{\mathbf{s}}}\right)
\right)
\tag{1}
\\&=
\sigma^{\mathbf{F}_{\Sigma^{\boldsymbol{\mathcal{A}}}}
\left(\mathbf{Pth}_{\boldsymbol{\mathcal{A}}}
\right)
}
\left(\left(
\mathrm{ip}^{(1,X)@}_{s_{j}}\left(
\mathrm{CH}^{(1)}_{s_{j}}\left(
\mathfrak{P}_{j}
\right)\right)
\right)_{j\in\bb{\mathbf{s}}}\right)
\tag{2}
\\&=
\sigma^{
\mathbf{Pth}_{\boldsymbol{\mathcal{A}}}
}
\left(\left(
\mathrm{ip}^{(1,X)@}_{s_{j}}\left(
\mathrm{CH}^{(1)}_{s_{j}}\left(
\mathfrak{P}_{j}
\right)\right)
\right)_{j\in\bb{\mathbf{s}}}\right).
\tag{3}
\end{align*}

In the just stated chain of equalities, the first equality recovers the value of the Curry-Howard mapping at $\mathfrak{P}$; the second equality holds since, according to Definition~\ref{DIp}, $\mathrm{ip}^{(1,X)@}$ is a $\Sigma^{\boldsymbol{\mathcal{A}}}$-homomorphism; finally, for the last equality let us note that, for every $j\in\bb{\mathbf{s}}$, $(\mathfrak{P}_{j},s_{j})$ is strictly smaller than $(\mathfrak{P},s)$ with respect to $\prec_{\mathbf{Pth}_{\boldsymbol{\mathcal{A}}}}$. Then, by the inductive hypothesis, we have that $\mathrm{ip}^{(1,X)@}_{s_{j}}(
\mathrm{CH}^{(1)}_{s_{j}}(
\mathfrak{P}_{j}
))$ is a path in $[\mathfrak{P}_{j}]_{s_{j}}$. Hence, the interpretation of the operation symbol $\sigma$ in $\mathbf{F}_{\Sigma^{\boldsymbol{\mathcal{A}}}}(
\mathbf{Pth}_{\boldsymbol{\mathcal{A}}}
)
$ becomes that of $\mathbf{Pth}_{\boldsymbol{\mathcal{A}}}$. Note that this operation is well-defined because paths in the same class have the same $(0,1)$-source and $(0,1)$-target according to Lemma~\ref{LCH}.

Hence $\mathrm{ip}^{(1,X)@}_{s}(\mathrm{CH}^{(1)}_{s}(\mathfrak{P}))$ is a path.

Since $\mathfrak{P}$ is a path of length at least one,  there exists at least one index $j\in\bb{\mathbf{s}}$ for which $\mathfrak{P}_{j}$ has length at least one. Then, since $(\mathfrak{P}_{j},\mathrm{ip}^{(1,X)@}_{s_{j}}(\mathrm{CH}^{(1)}_{s_{j}}(\mathfrak{P}_{j})))$ is a pair in $\mathrm{Ker}(\mathrm{CH}^{(1)})_{s_{j}}$, we have, by Lemma~\ref{LCH}, that $\mathrm{ip}^{(1,X)@}_{s_{j}}(\mathrm{CH}^{(1)}_{s_{j}}(\mathfrak{P}_{j}))$ is a path of length at least one. Then, since 
$$
\mathrm{ip}^{(1,X)@}_{s}\left(
\mathrm{CH}^{(1)}_{s}\left(
\mathfrak{P}
\right)\right)=
\sigma^{
\mathbf{Pth}_{\boldsymbol{\mathcal{A}}}
}
\left(\left(
\mathrm{ip}^{(1,X)@}_{s_{j}}\left(
\mathrm{CH}^{(1)}_{s_{j}}\left(
\mathfrak{P}_{j}
\right)\right)
\right)_{j\in\bb{\mathbf{s}}}\right),$$
we have, by Corollary~\ref{CPthWB}, that $\mathrm{ip}^{(1,X)@}_{s}(
\mathrm{CH}^{(1)}_{s}(
\mathfrak{P}
))$ is an echelonless path. Moreover, in virtue of Proposition~\ref{PRecov}, the  path extraction procedure from Lemma~\ref{LPthExtract} applied to $\mathrm{ip}^{(1,X)@}_{s}(
\mathrm{CH}^{(1)}_{s}(
\mathfrak{P}
))$ retrieves the family of paths $(
\mathrm{ip}^{(1,X)@}_{s_{j}}(
\mathrm{CH}^{(1)}_{s_{j}}(
\mathfrak{P}_{j}
))
)_{j\in\bb{\mathbf{s}}}$.

Hence, according to Definition~\ref{DCH}, the value of the Curry-Howard mapping at $\mathrm{ip}^{(1,X)@}_{s}(
\mathrm{CH}^{(1)}_{s}(
\mathfrak{P}
))$ is given by
\begin{flushleft}
$\mathrm{CH}^{(1)}_{s}\left(
\mathrm{ip}^{(1,X)@}_{s}\left(
\mathrm{CH}^{(1)}_{s}\left(
\mathfrak{P}
\right)\right)\right)$
\begin{align*}
\qquad
&=
\sigma^{\mathbf{T}_{\Sigma^{\boldsymbol{\mathcal{A}}}}(X)}
\left(\left(\mathrm{CH}^{(1)}_{s_{j}}\left(
\mathrm{ip}^{(1,X)@}_{s_{j}}\left(
\mathrm{CH}^{(1)}_{s_{j}}\left(
\mathfrak{P}_{j}
\right)\right)\right)
\right)_{j\in\bb{\mathbf{s}}}
\right)
\tag{1}
\\&=
\sigma^{\mathbf{T}_{\Sigma^{\boldsymbol{\mathcal{A}}}}(X)}
\left(\left(\mathrm{CH}^{(1)}_{s_{j}}\left(
\mathfrak{P}_{j}
\right)
\right)_{j\in\bb{\mathbf{s}}}
\right)
\tag{2}
\\&=
\mathrm{CH}^{(1)}_{s}\left(
\mathfrak{P}
\right).
\tag{3}
\end{align*}
\end{flushleft}

The first equality follows from the fact that $\mathrm{ip}^{(1,X)@}_{s}(
\mathrm{CH}^{(1)}_{s}(
\mathfrak{P}
))$ is an echelonless  path and Definition~\ref{DCH}; the second equality follows from the fact that, from the previous discussion, for every $j\in\bb{\mathbf{s}}$, $\mathrm{ip}^{(1,X)@}_{s_{j}}(\mathrm{CH}^{(1)}_{s_{j}}(\mathfrak{P}_{j}))$ is a path in $[\mathfrak{P}_{j}]_{s_{j}}$; finally, the last equality recovers the value of the  Curry-Howard mapping at $\mathfrak{P}$.

Therefore,
\begin{itemize}
\item[(i)] $\mathrm{ip}^{(1,X)@}_{s}(\mathrm{CH}^{(1)}_{s}(\mathfrak{P}))$ is a path in $\mathrm{Pth}_{\boldsymbol{\mathcal{A}},s}$; and
\item[(ii)] $\mathrm{ip}^{(1,X)@}_{s}(\mathrm{CH}^{(1)}_{s}(\mathfrak{P}))$ belongs to $[\mathfrak{P}]_{s}$.
\end{itemize}

The case of $\mathfrak{P}$ being an echelonless path follows.

This finishes the proof.
\end{proof}

\begin{remark} For a sort $s\in S$ and a path $\mathfrak{P}\in\mathrm{Pth}_{\boldsymbol{\mathcal{A}},s}$, 
$\mathrm{CH}^{(1)}_{s}(\mathfrak{P})$ is a syntactical description of $\mathfrak{P}$ and the path 
$\mathrm{ip}^{(1,X)@}_{s}(\mathrm{CH}^{(1)}_{s}(\mathfrak{P}))$ is the \emph{normalized path} arising from $\mathfrak{P}$, in which the rewrite rules, acting in parallel, are performed in the leftmost innermost possible position. This normalization process has guided us all the way up to here, see, e.g., Remark~\ref{RCH}.
\end{remark}

The following corollary states that $(1,0)$-identity paths are always normalized, or, what is equivalent, that $\mathrm{ip}^{(1,X)@}\circ\mathrm{CH}^{(1)}$ leaves invariant the $(1,0)$-identity paths.  

\begin{restatable}{corollary}{CIpId}
\label{CIpId}
Let $s$ be a sort in $S$, and $\mathfrak{P}$ a path in $\mathrm{Pth}_{\boldsymbol{\mathcal{A}},s}$. If $\mathfrak{P}$ is a  $(1,0)$-identity path, then $\mathrm{ip}^{(1,X)@}_{s}(\mathrm{CH}^{(1)}_{s}(\mathfrak{P}))=\mathfrak{P}$.
\end{restatable}
\begin{proof}
It follows from Proposition~\ref{PIpCH} and Corollary~\ref{CCHUZId}.
\end{proof}

We present a series of results that will be useful later on. The first result states that the composition 
$\mathrm{ip}^{(1,X)@}\circ \mathrm{CH}^{(1)}$ is a $\Sigma$-homomorphism.

\begin{restatable}{lemma}{LIpCHSigma}
\label{LIpCHSigma} Let $(\mathbf{s}, s)$ be a pair in  $S^{\star}\times S$, $\sigma$ an operation symbol in $\Sigma_{\mathbf{s},s}$ and $(\mathfrak{P}_{j})_{j\in\bb{\mathbf{s}}}$ a family of paths in $\mathrm{Pth}_{\boldsymbol{\mathcal{A}},\mathbf{s}}$.
Then the following equality holds
\allowdisplaybreaks
\begin{multline*}
\mathrm{ip}^{(1,X)@}_{s}\left(
\mathrm{CH}^{(1)}_{s}\left(
\sigma^{\mathbf{Pth}_{\boldsymbol{\mathcal{A}}}}\left(
\left(
\mathfrak{P}_{j}
\right)_{j\in\bb{\mathbf{s}}}
\right)\right)\right)\\=
\sigma^{\mathbf{Pth}_{\boldsymbol{\mathcal{A}}}}\left(
\left(
\mathrm{ip}^{(1,X)@}_{s_{j}}\left(
\mathrm{CH}^{(1)}_{s_{j}}\left(
\mathfrak{P}_{j}
\right)\right)
\right)_{j\in\bb{\mathbf{s}}}
\right).
\end{multline*}
\end{restatable}

\begin{proof}
It follows from the fact that $\mathrm{CH}^{(1)}$ and $\mathrm{ip}^{(1,X)@}$ are $\Sigma$-homomorphisms according to, respectively, Proposition~\ref{PCHHom} and Definition~\ref{DIp}.
\end{proof}

We next prove that the normalized path of a rewrite rule is the rewrite rule itself, or, what is equivalent, that $\mathrm{ip}^{(1,X)@}\circ \mathrm{CH}^{(1)}$ leaves invariant the rewrite rules.

\begin{restatable}{lemma}{LIpCHEch}
\label{LIpCHEch} Let $s$ be a sort in $S$ and $\mathfrak{p}$ a rewrite rule in $\mathcal{A}_{s}$. Then the following equality holds
\[
\mathrm{ip}^{(1,X)@}_{s}\left(
\mathrm{CH}^{(1)}_{s}\left(
\mathfrak{p}^{\mathbf{Pth}_{\boldsymbol{\mathcal{A}}}}
\right)\right)
=
\mathfrak{p}^{\mathbf{Pth}_{\boldsymbol{\mathcal{A}}}}.
\]
\end{restatable}
\begin{proof}
The following chain of equalities holds.
\allowdisplaybreaks
\begin{align*}
\mathrm{ip}^{(1,X)@}_{s}\left(
\mathrm{CH}^{(1)}_{s}\left(
\mathfrak{p}^{\mathbf{Pth}_{\boldsymbol{\mathcal{A}}}}
\right)\right)&=
\mathrm{ip}^{(1,X)@}_{s}\left(
\mathrm{CH}^{(1)}_{s}\left(
\mathrm{ech}^{(1,\mathcal{A})}_{s}\left(
\mathfrak{p}
\right)
\right)\right)
\tag{1}
\\&=
\mathrm{ip}^{(1,X)@}_{s}\left(
\eta^{(1,\mathcal{A})}_{s}\left(
\mathfrak{p}
\right)\right)
\tag{2}
\\&=
\mathrm{ech}^{(1,\mathcal{A})}_{s}\left(
\mathfrak{p}
\right)
\tag{3}
\\&=
\mathfrak{p}^{\mathbf{Pth}_{\boldsymbol{\mathcal{A}}}}.
\tag{4}
\end{align*}

In the just stated chain of equalities, the first equality unravels the interpretation of the constant operation symbol $\mathfrak{p}$ in the many-sorted partial $\Sigma^{\boldsymbol{\mathcal{A}}}$-algebra $\mathbf{Pth}_{\boldsymbol{\mathcal{A}}}$, according to Proposition~\ref{PPthCatAlg}; the second equality follows from Proposition~\ref{PCHA}; the third equality follows from Proposition~\ref{PIpEch}; finally, the last equality recovers the interpretation of the constant operation symbol $\mathfrak{p}$ in the many-sorted partial $\Sigma^{\boldsymbol{\mathcal{A}}}$-algebra $\mathbf{Pth}_{\boldsymbol{\mathcal{A}}}$, according to Proposition~\ref{PPthCatAlg}.

This completes the proof.
\end{proof}

We next prove that the normalized path of the $0$-source of a path is the $0$-source of the normalized path, or, what is equivalent, that $\mathrm{ip}^{(1,X)@}\circ \mathrm{CH}^{(1)}$ commutes with 
$\mathrm{sc}^{0\mathbf{Pth}_{\boldsymbol{\mathcal{A}}}}$. The same holds for the $0$-target.

\begin{restatable}{lemma}{LIpCHScTg}
\label{LIpCHScTg} Let $s$ be a sort in $S$ and let $\mathfrak{P}$ be a  path in $\mathrm{Pth}_{\boldsymbol{\mathcal{A}},s}$
then the following equality holds
\allowdisplaybreaks
\begin{align*}
\mathrm{ip}^{(1,X)@}_{s}\left(
\mathrm{CH}^{(1)}_{s}\left(
\mathrm{sc}^{0\mathbf{Pth}_{\boldsymbol{\mathcal{A}}}}_{s}\left(
\mathfrak{P}
\right)\right)\right)
&=
\mathrm{sc}^{0\mathbf{Pth}_{\boldsymbol{\mathcal{A}}}}_{s}\left(
\mathrm{ip}^{(1,X)@}_{s}\left(
\mathrm{CH}^{(1)}_{s}\left(
\mathfrak{P}
\right)\right)\right);
\\
\mathrm{ip}^{(1,X)@}_{s}\left(
\mathrm{CH}^{(1)}_{s}\left(
\mathrm{tg}^{0\mathbf{Pth}_{\boldsymbol{\mathcal{A}}}}_{s}\left(
\mathfrak{P}
\right)\right)\right)
&=
\mathrm{tg}^{0\mathbf{Pth}_{\boldsymbol{\mathcal{A}}}}_{s}\left(
\mathrm{ip}^{(1,X)@}_{s}\left(
\mathrm{CH}^{(1)}_{s}\left(
\mathfrak{P}
\right)\right)\right).
\end{align*}
\end{restatable}

\begin{proof}
We present the proof only for the first statement. The other one can be handled similarly.

The following chain of equalities holds
\allowdisplaybreaks
\begin{align*}
\mathrm{ip}^{(1,X)@}_{s}\left(
\mathrm{CH}^{(1)}_{s}\left(
\mathrm{sc}^{0\mathbf{Pth}_{\boldsymbol{\mathcal{A}}}}_{s}\left(
\mathfrak{P}
\right)\right)\right)&=
\mathrm{ip}^{(1,X)@}_{s}\left(
\mathrm{CH}^{(1)}_{s}\left(
\mathrm{ip}^{(1,0)\sharp}_{s}\left(
\mathrm{sc}^{(0,1)}_{s}\left(
\mathfrak{P}
\right)\right)\right)\right)
\tag{1}
\\&=
\mathrm{ip}^{(1,0)\sharp}_{s}\left(
\mathrm{sc}^{(0,1)}_{s}\left(
\mathfrak{P}
\right)
\right)
\tag{2}
\\&=
\mathrm{ip}^{(1,0)\sharp}_{s}\left(
\mathrm{sc}^{(0,1)}_{s}\left(
\mathrm{ip}^{(1,X)@}_{s}\left(
\mathrm{CH}^{(1)}_{s}\left(
\mathfrak{P}
\right)
\right)
\right)
\right)
\tag{3}
\\&=
\mathrm{sc}^{0\mathbf{Pth}_{\boldsymbol{\mathcal{A}}}}_{s}\left(
\mathrm{ip}^{(1,X)@}_{s}\left(
\mathrm{CH}^{(1)}_{s}\left(
\mathfrak{P}
\right)\right)
\right).
\tag{4}
\end{align*}

In the just stated chain of equalities, the first equality unravels the description of the operation symbol $\mathrm{sc}^{0}_{s}$ in the many-sorted partial $\Sigma^{\boldsymbol{\mathcal{A}}}$-algebra $\mathbf{Pth}_{\boldsymbol{\mathcal{A}}}$, according to Proposition~\ref{PPthCatAlg}; the second equality follows from Proposition~\ref{PIpUZ}; the third equality follows from Proposition~\ref{PIpCH} and Lemma~\ref{LCH}; finally, the last equality recovers the description of the operation symbol $\mathrm{sc}^{0}_{s}$ in the many-sorted partial $\Sigma^{\boldsymbol{\mathcal{A}}}$-algebra $\mathbf{Pth}_{\boldsymbol{\mathcal{A}}}$, according to Proposition~\ref{PPthCatAlg}.

This completes the proof.
\end{proof}

The following corollary states that one-step paths are always normalized, or, what is equivalent, that 
$\mathrm{ip}^{(1,X)@}\circ \mathrm{CH}^{(1)}$ leaves invariant the one-step paths. In this case there is just one derivation and it has been done, necessarily, in the leftmost innermost possible position.

\begin{restatable}{corollary}{CIpCHOneStep}
\label{CIpCHOneStep}
Let $s$ be a sort in $S$, and $\mathfrak{P}$ a path in $\mathrm{Pth}_{\boldsymbol{\mathcal{A}},s}$. If $\mathfrak{P}$ is a one-step path, then $\mathrm{ip}^{(1,X)@}_{s}(\mathrm{CH}^{(1)}_{s}(\mathfrak{P}))=\mathfrak{P}$.
\end{restatable}
\begin{proof}
It follows from Propositions~\ref{PIpCH} and~\ref{PCHOneStep}.
\end{proof}

In the following example we compute the free-completion of the identity path mapping applied to the value of the Curry-Howard mapping at a path.

\begin{example}[Continuation of Examples~\ref{ERun} and~\ref{ERunII}]\label{ERunIII}
We consider the path $\mathfrak{P}$ presented in Example~\ref{ERun}.

Let us recall from Example~\ref{ERunII} that the value of the Curry-Howard mapping at $\mathfrak{P}$ is the term in $\mathrm{T}_{\Sigma^{\boldsymbol{\mathcal{A}}}}(X)_{0}$ given by
\begin{multline*}
\mathrm{CH}^{(1)}_{0}\left(\mathfrak{P}\right)=
\Big(
\left(
\left(
\left(
\mathfrak{p}_{0}
\circ^{0}_{0}
\mathfrak{p}_{7}
\right)
\circ^{0}_{0}
\beta\left(
\mathfrak{p}_{4},
\delta\left(\mathfrak{q}_{0},\mathfrak{q}_{1}\right)
\circ^{0}_{0}
\mathfrak{p}_{3},
a
\right)
\right)
\circ^{0}_{0}
\mathfrak{p}_{5}
\right)
\circ_{0}^{0}
\\
\delta\left(x,\mathfrak{q}_{0}\circ_{1}^{0}\mathfrak{q}_{1}\right)
\circ^{0}_{0}\mathfrak{p}_{6}
\Big)
\circ^{0}_{0}
\beta\left(\mathfrak{p}_{4},\mathfrak{p}_{1}\circ^{0}_{0}\mathfrak{p}_{2},\mathfrak{p}_{1}\right).
\end{multline*}

Now we can apply $\mathrm{ip}^{(1,X)@}$ to $\mathrm{CH}^{(1)}_{0}(\mathfrak{P})$, to obtain a canonical path which is congruent to $\mathfrak{P}$ with respect to $\mathrm{Ker}(\mathrm{CH}^{(1)})$, and in which all rewrite rules are performed in the leftmost innermost possible position. This is, we recall, the normalized path that we have discussed in Proposition~\ref{PIpCH}. In what follows, for simplicity of notation, we let $\mathfrak{P}^{@}$ stand for 
$\mathrm{ip}^{(1,X)@}_{0}(\mathrm{CH}^{(1)}_{0}(\mathfrak{P}))$. The path $\mathfrak{P}^{@}$ is given by the following sequence. To simplify the presentation we have omitted  the superscript $\mathbf{T}_{\Sigma^{\boldsymbol{\mathcal{A}}}}(X)$.
\allowdisplaybreaks
\begin{align*}
\mathfrak{P}^{@}\colon
\beta(\gamma(z, a),b,b)&
\xymatrix@C=100pt{\ar[r]^-{\text{\Small{$\left(\mathfrak{p}_{4}, 
\beta(\underline{\quad}, b,b)
\right)$}}}&}
\beta(a, b, b) 
\\&
\xymatrix@C=100pt{\ar[r]^-{\text{\Small{$\left(\mathfrak{p}_{2}, 
\beta(a,\underline{\quad},b)
\right)$}}}&}
\beta(a, b, b) 
\\&
\xymatrix@C=100pt{\ar[r]^-{\text{\Small{$\left(\mathfrak{p}_{1}, 
\beta(a,\underline{\quad},b)
\right)$}}}&}
\beta(a, a, b) 
\\&
\xymatrix@C=100pt{\ar[r]^-{\text{\Small{$\left(\mathfrak{p}_{1}, 
\beta(a,a,\underline{\quad})
\right)$}}}&}
\beta(a,a,a) 
\\&
\xymatrix@C=100pt{\ar[r]^-{\text{\Small{$(\mathfrak{p}_{6}, 
\underline{\quad}
)$}}}&}
\delta(x,y) 
\\&
\xymatrix@C=100pt{\ar[r]^-{\text{\Small{$(\mathfrak{q}_{1}, 
\delta(x,\underline{\quad})
)$}}}&}
\delta(x,x) 
\\&
\xymatrix@C=100pt{\ar[r]^-{\text{\Small{$(\mathfrak{q}_{0}, 
\delta(x,\underline{\quad})
)$}}}&}
\delta(x,z) 
\\&
\xymatrix@C=100pt{\ar[r]^-{\text{\Small{$(\mathfrak{p}_{5}, 
\underline{\quad}
)$}}}&}
\beta(\gamma(z, a), b, a) 
\\&
\xymatrix@C=100pt{\ar[r]^-{\text{\Small{$(\mathfrak{p}_{4}, 
\beta(\underline{\quad}, b, a)
)$}}}&}
\beta(a, b, a) 
\\&
\xymatrix@C=100pt{\ar[r]^-{\text{\Small{$(\mathfrak{p}_{3}, 
\beta(a, \underline{\quad}, a)
)$}}}&}
\beta(a, \delta(x,y), a) 
\\&
\xymatrix@C=100pt{\ar[r]^-{\text{\Small{$(\mathfrak{q}_{0}, 
\beta(a, \delta(\underline{\quad},y), a)
)$}}}&}
\beta(a, \delta(z,y), a) 
\\&
\xymatrix@C=100pt{\ar[r]^-{\text{\Small{$(\mathfrak{q}_{1}, 
\beta(a, \delta(z,\underline{\quad}), a)
)$}}}&}
\beta(a, \delta(z,x), a)
\\&
\xymatrix@C=100pt{\ar[r]^-{\text{\Small{$(\mathfrak{p}_{7}, 
\underline{\quad}
)$}}}&}
\alpha
\\&
\xymatrix@C=100pt{\ar[r]^-{\text{\Small{$(\mathfrak{p}_{0}, 
\underline{\quad}
)$}}}&}
\delta(x,z).
\end{align*}

The path $\mathfrak{P}^{@}$ is a $(00000110001100)$-path in $\mathrm{Pth}_{\boldsymbol{\mathcal{A}},0}$ from $\beta(\gamma(z, a),b,b)$ to $\delta(x,z)$.

Let us note that $\mathfrak{P}^{@}$ has the same length, $(0,1)$-source, and $(0,1)$-target that $\mathfrak{P}$; that $\mathfrak{P}^{@}$ has the same intervals of echelons and echelonless subpaths  (and in the same positions) that $\mathfrak{P}$; and that, by Proposition~\ref{PIpCH}, $\mathrm{CH}^{(1)}_{0}(\mathfrak{P}) = \mathrm{CH}^{(1)}_{0}(\mathfrak{P}^{@})$.
\end{example}

\section{Order issues}

In this section we prove that $\coprod(\mathrm{ip}^{(1,X)@}\circ\mathrm{CH}^{(1)})$ is an order-preserving mapping from $\coprod\mathrm{Pth}_{\boldsymbol{\mathcal{A}}}$ to itself. But, to do so, we first prove that if the value of $\mathrm{CH}^{(1)}$ at a path $\leq_{\mathbf{T}_{\Sigma^{\boldsymbol{\mathcal{A}}}}(X)}$-precedes to the value of $\mathrm{CH}^{(1)}$ at another path, then the respective normalised paths satisfy the same inequality with respect to the order $\leq_{\mathbf{Pth}_{\boldsymbol{\mathcal{A}}}}$ introduced in Definition~\ref{DOrd}.

\begin{restatable}{proposition}{PIpCHOrd}
\label{PIpCHOrd} Let $s,t$ be sorts in $S$ and $\mathfrak{P}\in\mathrm{Pth}_{\boldsymbol{\mathcal{A}},s}$ and $\mathfrak{Q}\in\mathrm{Pth}_{\boldsymbol{\mathcal{A}},t}$. If 
\[
\left(
\mathrm{CH}^{(1)}_{t}\left(
\mathfrak{Q}
\right)
,t
\right)
\leq_{\mathbf{T}_{\Sigma^{\boldsymbol{\mathcal{A}}}}(X)}
\left(
\mathrm{CH}^{(1)}_{s}\left(
\mathfrak{P}
\right)
,s
\right),
\]
then the following inequality holds
\[
\left(
\mathrm{ip}^{(1,X)@}_{t}\left(
\mathrm{CH}^{(1)}_{t}\left(
\mathfrak{Q}\right)\right)
,t
\right)
\leq_{\mathbf{Pth}_{\boldsymbol{\mathcal{A}}}}
\left(
\mathrm{ip}^{(1,X)@}_{s}\left(
\mathrm{CH}^{(1)}_{s}\left(
\mathfrak{P}\right)\right)
,s
\right).
\]
\end{restatable}

\begin{proof}
Note that in case $t=s$ and $\mathrm{CH}^{(1)}_{s}(\mathfrak{Q})=\mathrm{CH}^{(1)}_{s}(\mathfrak{P})$, the statement follows by reflexivity of $\leq_{\mathbf{Pth}_{\boldsymbol{\mathcal{A}}}}$. Therefore, we can work under the the following assumption.
\begin{assumption}\label{AIpCHOrdI} The following strict inequality holds
\[
\left(
\mathrm{CH}^{(1)}_{t}\left(
\mathfrak{Q}
\right)
,t
\right)
<_{\mathbf{T}_{\Sigma^{\boldsymbol{\mathcal{A}}}}(X)}
\left(
\mathrm{CH}^{(1)}_{s}\left(
\mathfrak{P}
\right)
,s
\right).
\]
\end{assumption}

Moreover, the case in which $\mathfrak{P}$ is a $(1,0)$-identity path also follows immediately. Note that in this case 
$\mathfrak{P}=\mathrm{ip}^{(1,0)\sharp}_{s}(P)$
for some term $P\in \mathrm{T}_{\Sigma}(X)_{s}$. Hence, according to Proposition~\ref{PCHId}, we have that 
\[
\mathrm{CH}^{(1)}_{s}\left(\mathfrak{P}\right)=\eta^{(1,0)\sharp}_{s}\left(P\right).
\]

If we have that $\mathrm{CH}^{(1)}_{t}(\mathfrak{Q})$ is a subterm of sort $t$ of $\mathrm{CH}^{(1)}_{s}(\mathfrak{P})$, then we have that $\mathrm{CH}^{(1)}_{t}(\mathfrak{Q})$ is a term in $\eta^{(1,0)\sharp}[\mathrm{T}_{\Sigma}(X)]_{t}$. Hence, according to Corollary~\ref{CCHId}, $\mathfrak{Q}$ is also a $(1,0)$-identity path, i.e., $\mathfrak{Q}=\mathrm{ip}^{(1,0)\sharp}_{t}(Q)$
for some term $Q\in \mathrm{T}_{\Sigma}(X)_{t}$. Note also that, by Proposition~\ref{PCHId}, we have that 
\[
\mathrm{CH}^{(1)}_{t}\left(\mathfrak{Q}\right)=\eta^{(1,0)\sharp}_{t}\left(Q\right).
\]

Since $
(
\mathrm{CH}^{(1)}_{t}(
\mathfrak{Q}
)
,t
)
\leq_{\mathbf{T}_{\Sigma^{\boldsymbol{\mathcal{A}}}}(X)}
(
\mathrm{CH}^{(1)}_{s}(
\mathfrak{P}
)
,s
)$, we conclude that \[
\left(
\eta^{(1,0)\sharp}_{t}(Q)
,t
\right)
\leq_{\mathbf{T}_{\Sigma^{\boldsymbol{\mathcal{A}}}}(X)}
\left(
\eta^{(1,0)\sharp}_{s}(P)
,s
\right).\]
Thus, $(Q,t)\leq_{\mathbf{T}_{\Sigma}(X)} (P,s)$.  Following Definition~\ref{DOrd}, we conclude that 
\[
\left(
\mathfrak{Q},t
\right)
\leq_{\mathbf{Pth}_{\boldsymbol{\mathcal{A}}}}
\left(
\mathfrak{P},s
\right).
\]

Finally note that, according to Proposition~\ref{PIpUZ}, we have that 
\allowdisplaybreaks
\begin{align*}
\mathrm{ip}^{(1,X)@}_{t}\left(
\mathrm{CH}^{(1)}_{t}\left(
\mathfrak{Q}\right)\right)
&=
\mathfrak{Q};
&
\mathrm{ip}^{(1,X)@}_{s}\left(
\mathrm{CH}^{(1)}_{s}\left(
\mathfrak{P}\right)\right)
&=
\mathfrak{P}.
\end{align*}

This case follows. 

Hence we can work under the following assumption.
\begin{assumption}\label{AIpCHOrdII}
$\mathfrak{P}$ is not a $(1,0)$-identity path.
\end{assumption}

We prove the statement by Artinian induction on $(\coprod\mathrm{Pth}_{\boldsymbol{\mathcal{A}}}, \leq_{\mathbf{Pth}_{\boldsymbol{\mathcal{A}}}})$.

\textsf{Base step of the Artinian induction.}

Let $(\mathfrak{P},s)$ be a minimal element of $(\coprod\mathrm{Pth}_{\boldsymbol{\mathcal{A}}}, \leq_{\mathbf{Pth}_{\boldsymbol{\mathcal{A}}}})$. Then, by Proposition~\ref{PMinimal} and taking into account that $\mathfrak{P}$ is a non-$(1,0)$-identity path, we have that $\mathfrak{P}$ is an echelon associated to a rewriting rule $\mathfrak{p}\in\mathcal{A}_{s}$, i.e., $\mathfrak{P}=\mathrm{ech}^{(1,\mathcal{A})}_{s}(\mathfrak{p})$.

In this case, we have according to Definition~\ref{DCH}, that 
\[
\mathrm{CH}^{(1)}_{s}\left(\mathfrak{P}\right)=\mathfrak{p}^{\mathbf{T}_{\Sigma^{\boldsymbol{\mathcal{A}}}}(X)}.
\]

Note that $\mathfrak{p}^{\mathbf{T}_{\Sigma^{\boldsymbol{\mathcal{A}}}}(X)}$ is a term in $\mathrm{B}^{0}_{\Sigma^{\boldsymbol{\mathcal{A}}}}(X)_{s}$, therefore the only possibility for $\mathrm{CH}^{(1)}_{t}(\mathfrak{Q})$ of being a subterm of sort $t$ of $\mathfrak{p}^{\mathbf{T}_{\Sigma^{\boldsymbol{\mathcal{A}}}}(X)}$ is that $t=s$ and $\mathrm{CH}^{(1)}_{s}(\mathfrak{Q})=\mathfrak{p}^{\mathbf{T}_{\Sigma^{\boldsymbol{\mathcal{A}}}}(X)}$. This case follows by reflexivity of $\leq_{\mathbf{Pth}_{\boldsymbol{\mathcal{A}}}}$.

This completes the base step.

\textsf{Inductive step of the Artinian induction.}

Let $(\mathfrak{P},s)$ be a non-minimal element in $(\coprod\mathrm{Pth}_{\boldsymbol{\mathcal{A}}},\leq_{\mathbf{Pth}_{\boldsymbol{\mathcal{A}}}})$. Let us suppose that, for every sort $u\in S$ and every path $\mathfrak{R}\in\mathrm{Pth}_{\boldsymbol{\mathcal{A}},u}$, if $(\mathfrak{R},u)<_{\mathbf{Pth}_{\boldsymbol{\mathcal{A}}}} (\mathfrak{P},s)$, then the statement holds for $\mathfrak{R}$, i.e., if
the following inequality holds
\[
\left(
\mathrm{CH}^{(1)}_{t}\left(
\mathfrak{Q}
\right)
,t
\right)
\leq_{\mathbf{T}_{\Sigma^{\boldsymbol{\mathcal{A}}}}(X)}
\left(
\mathrm{CH}^{(1)}_{u}\left(
\mathfrak{R}
\right)
,u
\right),
\]
then the following inequality holds
\[
\left(
\mathrm{ip}^{(1,X)@}_{t}\left(
\mathrm{CH}^{(1)}_{t}\left(
\mathfrak{Q}\right)\right)
,t
\right)
\leq_{\mathbf{Pth}_{\boldsymbol{\mathcal{A}}}}
\left(
\mathrm{ip}^{(1,X)@}_{u}\left(
\mathrm{CH}^{(1)}_{u}\left(
\mathfrak{R}\right)\right)
,u
\right).
\]

Since $(\mathfrak{P},s)$ is a non-minimal element of $(\coprod\mathrm{Pth}_{\boldsymbol{\mathcal{A}}}, \leq_{\mathbf{Pth}_{\boldsymbol{\mathcal{A}}}})$ and taking into account that, under Assumption~\ref{AIpCHOrdII}, we are assuming that $\mathfrak{P}$ is a non-$(1,0)$-identity path, we have, by Lemma~\ref{LOrdI}, that $\mathfrak{P}$ is either (1) a path of length strictly greater than one containing at least one echelon or (2) and echelonless path.

If (1), then let $i\in\bb{\mathfrak{P}}$ be the first index for which the one-step subpath $\mathfrak{P}^{i,i}$ is an echelon. We distinguish two cases accordingly. 

If (1.1) $i=0$, i.e., if $\mathfrak{P}$ has its first echelon on its first step, then according to Definition~\ref{DCH}, we have that 
\[
\mathrm{CH}^{(1)}_{s}\left(
\mathfrak{P}
\right)
=
\mathrm{CH}^{(1)}_{s}\left(
\mathfrak{P}^{1,\bb{\mathfrak{P}}-1}
\right)
\circ^{0\mathbf{T}_{\Sigma^{\boldsymbol{\mathcal{A}}}}(X)}_{s}
\mathrm{CH}^{(1)}_{s}\left(
\mathfrak{P}^{0,0}
\right).
\]

Since $\mathrm{ip}^{(1,X)@}$ is a $\Sigma^{\boldsymbol{\mathcal{A}}}$-homomorphism according to Definition~\ref{DIp}, the following equality holds
\[
\mathrm{ip}^{(1,X)@}_{s}\left(
\mathrm{CH}^{(1)}_{s}\left(
\mathfrak{P}
\right)
\right)=
\mathrm{ip}^{(1,X)@}_{s}\left(
\mathrm{CH}^{(1)}_{s}\left(
\mathfrak{P}^{1,\bb{\mathfrak{P}}-1}
\right)
\right)
\circ^{0\mathbf{Pth}_{\boldsymbol{\mathcal{A}}}}
\mathrm{ip}^{(1,X)@}_{s}\left(
\mathrm{CH}^{(1)}_{s}\left(
\mathfrak{P}^{0,0}
\right)
\right).\]

Now taking into account Proposition~\ref{PIpCH} and Lemma~\ref{LCHEchInt}, we have that $\mathrm{ip}^{(1,X)@}_{s}(
\mathrm{CH}^{(1)}_{s}(
\mathfrak{P}
)
)$ is a path of length strictly greater than one containing its first echelon on its first step. Moreover, taking into account Proposition~\ref{PIpCH} and Lemma~\ref{LCH}, we have that 
\allowdisplaybreaks
\begin{align*}
\left(\mathrm{ip}^{(1,X)@}_{s}\left(
\mathrm{CH}^{(1)}_{s}\left(
\mathfrak{P}
\right)
\right)\right)^{0,0}
&=
\mathrm{ip}^{(1,X)@}_{s}\left(
\mathrm{CH}^{(1)}_{s}\left(
\mathfrak{P}^{0,0}
\right)
\right);
\\
\left(\mathrm{ip}^{(1,X)@}_{s}\left(
\mathrm{CH}^{(1)}_{s}\left(
\mathfrak{P}
\right)
\right)\right)^{1,\bb{\mathfrak{P}}-1}
&=
\mathrm{ip}^{(1,X)@}_{s}\left(
\mathrm{CH}^{(1)}_{s}\left(
\mathfrak{P}^{1,\bb{\mathfrak{P}}-1}
\right)
\right).
\end{align*}

According to Definition~\ref{DOrd}, the following inequalities hold
\[
\left(
\mathrm{ip}^{(1,X)@}_{s}\left(
\mathrm{CH}^{(1)}_{s}\left(
\mathfrak{P}^{0,0}
\right)
\right),
s
\right)
\leq_{\mathbf{Pth}_{\boldsymbol{\mathcal{A}}}}
\left(
\mathrm{ip}^{(1,X)@}_{s}\left(
\mathrm{CH}^{(1)}_{s}\left(
\mathfrak{P}
\right)
\right),
s
\right)
\]
\[
\left(
\mathrm{ip}^{(1,X)@}_{s}\left(
\mathrm{CH}^{(1)}_{s}\left(
\mathfrak{P}^{1,\bb{\mathfrak{P}}-1}
\right)
\right),
s
\right)
\leq_{\mathbf{Pth}_{\boldsymbol{\mathcal{A}}}}
\left(
\mathrm{ip}^{(1,X)@}_{s}\left(
\mathrm{CH}^{(1)}_{s}\left(
\mathfrak{P}
\right)
\right),
s
\right).
\]

Now, if $\mathfrak{Q}$ is a path in $\mathrm{Pth}_{\boldsymbol{\mathcal{A}},t}$ satisfying that 
\[
\left(
\mathrm{CH}^{(1)}_{t}\left(
\mathfrak{Q}
\right)
,t
\right)
\leq_{\mathbf{T}_{\Sigma^{\boldsymbol{\mathcal{A}}}}(X)}
\left(
\mathrm{CH}^{(1)}_{s}\left(
\mathfrak{P}
\right)
,s
\right),
\]
then, under Assumption~\ref{AIpCHOrdI}, either one of the following conditions must hold
\begin{itemize}
\item[(1.1.1)]  $t=s$ and $\mathrm{CH}^{(1)}_{t}(\mathfrak{Q})=\mathrm{CH}^{(1)}_{s}(\mathfrak{P}^{0,0})$;
\item[(1.1.2)] $t=s$ and $\mathrm{CH}^{(1)}_{t}(\mathfrak{Q})=\mathrm{CH}^{(1)}_{s}(\mathfrak{P}^{1,\bb{\mathfrak{P}}-1})$;
\item[(1.1.3)] $( \mathrm{CH}^{(1)}_{t}(\mathfrak{Q}), t)\leq_{\mathbf{T}_{\Sigma^{\boldsymbol{\mathcal{A}}}}(X)} (\mathrm{CH}^{(1)}_{s}(\mathfrak{P}^{0,0}), s)$;
\item[(1.1.4)] $( \mathrm{CH}^{(1)}_{t}(\mathfrak{Q}), t)\leq_{\mathbf{T}_{\Sigma^{\boldsymbol{\mathcal{A}}}}(X)} (\mathrm{CH}^{(1)}_{s}(\mathfrak{P}^{1,\bb{\mathfrak{P}}-1}), s)$.
\end{itemize}

Note that cases (1.1.1) and (1.1.2) have already been discussed. The cases (1.1.3) and (1.1.4) follow by induction and the transitivity of $\leq_{\mathbf{Pth}_{\boldsymbol{\mathcal{A}}}}$.

This completes Case~(1.1).

If (1.2) $i>0$, i.e., if $\mathfrak{P}$ has its first echelon on a step different from the initial one, then according to Definition~\ref{DCH}, we have that 
\[
\mathrm{CH}^{(1)}_{s}\left(
\mathfrak{P}
\right)
=
\mathrm{CH}^{(1)}_{s}\left(
\mathfrak{P}^{i,\bb{\mathfrak{P}}-1}
\right)
\circ^{0\mathbf{T}_{\Sigma^{\boldsymbol{\mathcal{A}}}}(X)}_{s}
\mathrm{CH}^{(1)}_{s}\left(
\mathfrak{P}^{0,i-1}
\right).
\]

Since $\mathrm{ip}^{(1,X)@}$ is a $\Sigma^{\boldsymbol{\mathcal{A}}}$-homomorphism according to Definition~\ref{DIp}, the following equality holds
\[
\mathrm{ip}^{(1,X)@}_{s}\left(
\mathrm{CH}^{(1)}_{s}\left(
\mathfrak{P}
\right)
\right)=
\mathrm{ip}^{(1,X)@}_{s}\left(
\mathrm{CH}^{(1)}_{s}\left(
\mathfrak{P}^{i,\bb{\mathfrak{P}}-1}
\right)
\right)
\circ^{0\mathbf{Pth}_{\boldsymbol{\mathcal{A}}}}
\mathrm{ip}^{(1,X)@}_{s}\left(
\mathrm{CH}^{(1)}_{s}\left(
\mathfrak{P}^{0,i-1}
\right)
\right).\]

Now taking into account Proposition~\ref{PIpCH} and Lemma~\ref{LCHEchNInt}, we have that $\mathrm{ip}^{(1,X)@}_{s}(
\mathrm{CH}^{(1)}_{s}(
\mathfrak{P}
)
)$ is a path of length strictly greater than one containing its first echelon on the step $i\in\bb{\mathfrak{P}}$, where $i>0$. Moreover, taking into account Proposition~\ref{PIpCH} and Lemma~\ref{LCH}, we have that 
\allowdisplaybreaks
\begin{align*}
\left(\mathrm{ip}^{(1,X)@}_{s}\left(
\mathrm{CH}^{(1)}_{s}\left(
\mathfrak{P}
\right)
\right)\right)^{0,i-1}
&=
\mathrm{ip}^{(1,X)@}_{s}\left(
\mathrm{CH}^{(1)}_{s}\left(
\mathfrak{P}^{0,i-1}
\right)
\right);
\\
\left(\mathrm{ip}^{(1,X)@}_{s}\left(
\mathrm{CH}^{(1)}_{s}\left(
\mathfrak{P}
\right)
\right)\right)^{i,\bb{\mathfrak{P}}-1}
&=
\mathrm{ip}^{(1,X)@}_{s}\left(
\mathrm{CH}^{(1)}_{s}\left(
\mathfrak{P}^{i,\bb{\mathfrak{P}}-1}
\right)
\right).
\end{align*}

According to Definition~\ref{DOrd}, the following inequalities hold
\[
\left(
\mathrm{ip}^{(1,X)@}_{s}\left(
\mathrm{CH}^{(1)}_{s}\left(
\mathfrak{P}^{0,i-1}
\right)
\right),
s
\right)
\leq_{\mathbf{Pth}_{\boldsymbol{\mathcal{A}}}}
\left(
\mathrm{ip}^{(1,X)@}_{s}\left(
\mathrm{CH}^{(1)}_{s}\left(
\mathfrak{P}
\right)
\right),
s
\right)
\]
\[
\left(
\mathrm{ip}^{(1,X)@}_{s}\left(
\mathrm{CH}^{(1)}_{s}\left(
\mathfrak{P}^{i,\bb{\mathfrak{P}}-1}
\right)
\right),
s
\right)
\leq_{\mathbf{Pth}_{\boldsymbol{\mathcal{A}}}}
\left(
\mathrm{ip}^{(1,X)@}_{s}\left(
\mathrm{CH}^{(1)}_{s}\left(
\mathfrak{P}
\right)
\right),
s
\right).
\]
\\

Now, if $\mathfrak{Q}$ is a path in $\mathrm{Pth}_{\boldsymbol{\mathcal{A}},t}$ satisfying that 
\[
\left(
\mathrm{CH}^{(1)}_{t}\left(
\mathfrak{Q}
\right)
,t
\right)
\leq_{\mathbf{T}_{\Sigma^{\boldsymbol{\mathcal{A}}}}(X)}
\left(
\mathrm{CH}^{(1)}_{s}\left(
\mathfrak{P}
\right)
,s
\right),
\]
then, under Assumption~\ref{AIpCHOrdI}, either one of the following conditions must hold
\begin{itemize}
\item[(1.2.1)]  $t=s$ and $\mathrm{CH}^{(1)}_{t}(\mathfrak{Q})=\mathrm{CH}^{(1)}_{s}(\mathfrak{P}^{0,i-1})$;
\item[(1.2.2)] $t=s$ and $\mathrm{CH}^{(1)}_{t}(\mathfrak{Q})=\mathrm{CH}^{(1)}_{s}(\mathfrak{P}^{i,\bb{\mathfrak{P}}-1})$;
\item[(1.2.3)] $( \mathrm{CH}^{(1)}_{t}(\mathfrak{Q}), t)\leq_{\mathbf{T}_{\Sigma^{\boldsymbol{\mathcal{A}}}}(X)} (\mathrm{CH}^{(1)}_{s}(\mathfrak{P}^{0,i-1}), s)$;
\item[(1.2.4)] $( \mathrm{CH}^{(1)}_{t}(\mathfrak{Q}), t)\leq_{\mathbf{T}_{\Sigma^{\boldsymbol{\mathcal{A}}}}(X)} (\mathrm{CH}^{(1)}_{s}(\mathfrak{P}^{i,\bb{\mathfrak{P}}-1}), s)$.
\end{itemize}

Note that cases (1.2.1) and (1.2.2) have already been discussed. The cases (1.2.3) and (1.2.4) follow by induction and the transitivity of $\leq_{\mathbf{Pth}_{\boldsymbol{\mathcal{A}}}}$.

This completes Case~(1.2).

If~(2), i.e., if $\mathfrak{P}$ is an echelonless path, then in virtue of Lemma~\ref{LPthHeadCt} there exists a unique word $\mathbf{s}\in S^{\star}-\{\lambda\}$ and a unique operation symbol $\sigma\in \Sigma_{\mathbf{s},s}$ associated to $\mathfrak{P}$. Let $(\mathfrak{P}_{j})_{j\in\bb{\mathbf{s}}}$ be the family of paths in $\mathrm{Pth}_{\boldsymbol{\mathcal{A}},\mathbf{s}}$ which, in virtue of Lemma~\ref{LPthExtract}, we can extract from $\mathfrak{P}$. Then according to Definition~\ref{DCH}, the value of the Curry-Howard mapping at $\mathfrak{P}$ is 
\[
\mathrm{CH}^{(1)}_{s}\left(\mathfrak{P}\right)=
\sigma^{\mathbf{T}_{\Sigma^{\boldsymbol{\mathcal{A}}}}(X)}\left(
\left(
\mathrm{CH}^{(1)}_{s_{j}}\left(
\mathfrak{P}_{j}
\right)
\right)_{j\in\bb{\mathbf{s}}}
\right).
\]

Since $\mathrm{ip}^{(1,X)@}$ is a $\Sigma^{\boldsymbol{\mathcal{A}}}$-homomorphism according to Definition~\ref{DIp}, the following equality holds
\[
\mathrm{ip}^{(1,X)@}_{s}\left(
\mathrm{CH}^{(1)}_{s}\left(
\mathfrak{P}
\right)
\right)=
\sigma^{\mathbf{Pth}_{\boldsymbol{\mathcal{A}}}}
\left(\left(
\mathrm{ip}^{(1,X)@}_{s_{j}}\left(
\mathrm{CH}^{(1)}_{s_{j}}\left(
\mathfrak{P}_{j}
\right)
\right)
\right)_{j\in\bb{\mathbf{s}}}
\right)
.\]

Now taking into account Proposition~\ref{PIpCH} and Lemma~\ref{LCHNEch}, we have that $\mathrm{ip}^{(1,X)@}_{s}(
\mathrm{CH}^{(1)}_{s}(
\mathfrak{P}
)
)$ is an echelonless path associated to the operation symbol $\sigma\in \Sigma_{\mathbf{s},s}$. Moreover, taking into account Proposition~\ref{PIpCH}, Lemma~\ref{LCH}, Corollary~\ref{CPthWB} and Proposition~\ref{PRecov}, we have that the path extraction procedure from Lemma~\ref{LPthExtract} applied to the path $\mathrm{ip}^{(1,X)@}_{s}(
\mathrm{CH}^{(1)}_{s}(
\mathfrak{P}
)
)$ retrieves the family of paths
\[\left(
\mathrm{ip}^{(1,X)@}_{s_{j}}\left(
\mathrm{CH}^{(1)}_{s_{j}}\left(
\mathfrak{P}_{j}
\right)
\right)
\right)_{j\in\bb{\mathbf{s}}}.\]

According to Definition~\ref{DOrd}, for every $j\in\bb{\mathbf{s}}$, the following inequality holds
\allowdisplaybreaks
\begin{align*}
\left(
\mathrm{ip}^{(1,X)@}_{s_{j}}\left(
\mathrm{CH}^{(1)}_{s_{j}}\left(
\mathfrak{P}_{j}
\right)
\right),
s_{j}
\right)
&\leq_{\mathbf{Pth}_{\boldsymbol{\mathcal{A}}}}
\left(
\mathrm{ip}^{(1,X)@}_{s}\left(
\mathrm{CH}^{(1)}_{s}\left(
\mathfrak{P}
\right)
\right),
s
\right).
\end{align*}

Now, if $\mathfrak{Q}$ is a path in $\mathrm{Pth}_{\boldsymbol{\mathcal{A}},t}$ satisfying that 
\[
\left(
\mathrm{CH}^{(1)}_{t}\left(
\mathfrak{Q}
\right)
,t
\right)
\leq_{\mathbf{T}_{\Sigma^{\boldsymbol{\mathcal{A}}}}(X)}
\left(
\mathrm{CH}^{(1)}_{s}\left(
\mathfrak{P}
\right)
,s
\right),
\]
then, under Assumption~\ref{AIpCHOrdI}, either one of the following conditions must hold
\begin{itemize}
\item[(2.1)]  there exists some $j\in\bb{\mathbf{s}}$ for which $t=s_{j}$ and $\mathrm{CH}^{(1)}_{t}(\mathfrak{Q})=\mathrm{CH}^{(1)}_{s_{j}}(\mathfrak{P}_{j})$;
\item[(2.2)] there exists some $j\in\bb{\mathbf{s}}$ for which  $( \mathrm{CH}^{(1)}_{t}(\mathfrak{Q}), t)\leq_{\mathbf{T}_{\Sigma^{\boldsymbol{\mathcal{A}}}}(X)} (\mathrm{CH}^{(1)}_{s_{j}}(\mathfrak{P}_{j}, s_{j})$.
\end{itemize}

Note that Case (2.1) has already been discussed. The Case (2.2) follows by induction and the transitivity of $\leq_{\mathbf{Pth}_{\boldsymbol{\mathcal{A}}}}$.

This completes Case~(2).

This completes the proof.
\end{proof}

The significance of the just stated proposition for our purposes is captured in the following corollary.

\begin{restatable}{corollary}{CIpCHOrd}
\label{CIpCHOrd} The mapping $\coprod(\mathrm{ip}^{(1,X)@}\circ \mathrm{CH}^{(1)})$ from $\coprod\mathrm{Pth}_{\boldsymbol{\mathcal{A}}}$ to $\coprod\mathrm{Pth}_{\boldsymbol{\mathcal{A}}}$ is order-preserving
\[
\textstyle
\coprod(\mathrm{ip}^{(1,X)@}\circ \mathrm{CH}^{(1)})
\colon
\left(
\coprod\mathrm{Pth}_{\boldsymbol{\mathcal{A}}},
\leq_{\mathbf{Pth}_{\boldsymbol{\mathcal{A}}}}
\right)
\mor
\left(
\coprod\mathrm{Pth}_{\boldsymbol{\mathcal{A}}},
\leq_{\mathbf{Pth}_{\boldsymbol{\mathcal{A}}}}
\right),
\]
that is, given pairs $(\mathfrak{Q},t)$, $(\mathfrak{P},s)$ in $\coprod\mathrm{Pth}_{\boldsymbol{\mathcal{A}}}$, if $(\mathfrak{Q},t)\leq_{\mathbf{Pth}_{\boldsymbol{\mathcal{A}}}}(\mathfrak{P},s)$ then
\[
\left(
\mathrm{ip}^{(1,X)@}_{t}\left(
\mathrm{CH}^{(1)}_{t}\left(
\mathfrak{Q}
\right)\right),t\right)
\leq_{\mathbf{Pth}_{\boldsymbol{\mathcal{A}}}}
\left(
\mathrm{ip}^{(1,X)@}_{s}\left(
\mathrm{CH}^{(1)}_{s}\left(
\mathfrak{P}
\right)\right),s\right).
\]
\end{restatable}
\begin{proof}
It follows from Proposition~\ref{PCHMono} and Proposition~\ref{PIpCHOrd}.
\end{proof}

\chapter{The correspondence between paths and terms}\label{S1H}

In this chapter we define a binary relation $\Theta^{(1)}$ on $\mathrm{T}_{\Sigma^{\boldsymbol{\mathcal{A}}}}(X)$  with the objective of matching different terms that, by $\mathrm{ip}^{(1,X)@}$, are sent to paths in the same equivalence class relative to $\mathrm{Ker}(\mathrm{CH}^{(1)})$. We also consider the smallest $\Sigma^{\boldsymbol{\mathcal{A}}}$-congruence containing $\Theta^{(1)}$, denoted by $\Theta^{[1]}$, and we check that this congruence also has the property of matching different terms that, by $\mathrm{ip}^{(1,X)@}$, are sent to paths in the same equivalence class relative to $\mathrm{Ker}(\mathrm{CH}^{(1)})$. We then consider the quotient $\Sigma^{\boldsymbol{\mathcal{A}}}$-algebra $\mathbf{T}_{\Sigma^{\boldsymbol{\mathcal{A}}}}(X)/{\Theta^{[1]}}$ and show that for the canonical projection 
$\mathrm{pr}^{\Theta^{[1]}}$ from $\mathbf{T}_{\Sigma^{\boldsymbol{\mathcal{A}}}}(X)$ to $\mathbf{T}_{\Sigma^{\boldsymbol{\mathcal{A}}}}(X)/{\Theta^{[1]}}$, the composition $\mathrm{pr}^{\Theta^{[1]}}\circ \mathrm{CH}^{(1)}$ is a $\Sigma^{\boldsymbol{\mathcal{A}}}$-homomorphism. Moreover, we show that if a term $P$ in $\mathrm{T}_{\Sigma^{\boldsymbol{\mathcal{A}}}}(X)_{s}$ is such that $\mathrm{ip}^{(1,X)@}_{s}(P)$ is a path in $\mathrm{Pth}_{\boldsymbol{\mathcal{A}},s}$, then the terms $P$ and $\mathrm{CH}^{(1)}_{s}(\mathrm{ip}^{(1,X)@}_{s}(P))$ are $\Theta^{[1]}_{s}$-related. Indeed, for a pair of terms $(P,Q)$ in $\Theta^{[1]}_{s}$, we show that $\mathrm{ip}^{(1,X)@}_{s}(P)$ is a path in $\mathrm{Pth}_{\boldsymbol{\mathcal{A}},s}$ if, and only if, $\mathrm{ip}^{(1,X)@}_{s}(Q)$ is a path in $\mathrm{Pth}_{\boldsymbol{\mathcal{A}},s}$. Moreover, when one of such conditions holds, the two paths $\mathrm{ip}^{(1,X)@}_{s}(P)$ and $\mathrm{ip}^{(1,X)@}_{s}(Q)$ have the same image under $\mathrm{CH}^{(1)}$. This, in particular entails that, if a term $P$ is $\Theta^{[1]}_{s}$-related with $\mathrm{CH}^{(1)}_{s}(\mathfrak{P})$, then $\mathrm{ip}^{(1,X)@}_{s}(P)$  is itself a path in the equivalence class $[\mathfrak{P}]_{s}$. Furthermore we show that, for two paths $\mathfrak{P}$ and 
$\mathfrak{P}'$ in $\mathrm{Pth}_{\boldsymbol{\mathcal{A}},s}$, if $\mathrm{CH}^{(1)}_{s}(\mathfrak{P})$ and 
$\mathrm{CH}^{(1)}_{s}(\mathfrak{P}')$ are $\Theta^{[1]}_{s}$-related then the equivalence classes 
$[\mathfrak{P}]_{s}$ and $[\mathfrak{P}']_{s}$ are equal.


We begin by defining a binary relation $\Theta^{(1)}$ on $\mathrm{T}_{\Sigma^{\boldsymbol{\mathcal{A}}}}(X)$ that intertwines $0$-sources, $0$-targets and $0$-compositions with $\mathrm{CH}^{(1)}$. This relation, as well as other relations associated with it, will play a major role in the subsequent work.
\begin{restatable}{definition}{DTheta}
\label{DTheta} 
\index{Theta!first-order!$\Theta^{(1)}$}
We let $\Theta^{(1)}=(\Theta^{(1)}_{s})_{s\in S}$ stand for the binary relation on $\mathrm{T}_{\Sigma^{\boldsymbol{\mathcal{A}}}}(X)$ consisting exactly of the following pairs of terms:
\begin{itemize}
\item[(i)] For every sort $s\in S$ and every path $\mathfrak{P}$ in $\mathrm{Pth}_{\boldsymbol{\mathcal{A}},s}$, 
$$
\left(\mathrm{CH}^{(1)}_{s}\left(
\mathrm{sc}^{0\mathbf{Pth}_{\boldsymbol{\mathcal{A}}}}_{s}\left(
\mathfrak{P}
\right)\right),
\mathrm{sc}^{0\mathbf{T}_{\Sigma^{\boldsymbol{\mathcal{A}}}}(X)}_{s}\left(
\mathrm{CH}^{(1)}_{s}\left(
\mathfrak{P}
\right)\right)
\right)\in\Theta^{(1)}_{s};
$$ 
\item[(ii)] For every sort $s\in S$ and every path $\mathfrak{P}$ in $\mathrm{Pth}_{\boldsymbol{\mathcal{A}},s}$, 
$$
\left(\mathrm{CH}^{(1)}_{s}\left(
\mathrm{tg}^{0\mathbf{Pth}_{\boldsymbol{\mathcal{A}}}}_{s}\left(
\mathfrak{P}
\right)\right),
\mathrm{tg}^{0\mathbf{T}_{\Sigma^{\boldsymbol{\mathcal{A}}}}(X)}_{s}\left(
\mathrm{CH}^{(1)}_{s}\left(
\mathfrak{P}
\right)\right)
\right)\in\Theta^{(1)}_{s};
$$ 
\item[(iii)] For every sort $s\in S$ and paths $\mathfrak{Q},\mathfrak{P}$ in $\mathrm{Pth}_{\boldsymbol{\mathcal{A}},s}$, if
$\mathrm{sc}^{(0,1)}_{s}(\mathfrak{Q})=\mathrm{tg}^{(0,1)}_{s}(\mathfrak{P})$, 
$$
\left(\mathrm{CH}^{(1)}_{s}\left(
\mathfrak{Q}
\circ^{0\mathbf{Pth}_{\boldsymbol{\mathcal{A}}}}_{s}
\mathfrak{P}
\right),
\mathrm{CH}^{(1)}_{s}\left(
\mathfrak{Q}
\right)
\circ^{0\mathbf{T}_{\Sigma^{\boldsymbol{\mathcal{A}}}}(X)}_{s}
\mathrm{CH}^{(1)}_{s}\left(
\mathfrak{P}
\right)
\right)\in\Theta^{(1)}_{s}.
$$ 
\end{itemize}
This completes the definition of $\Theta^{(1)}$.
\end{restatable}

Let us recall that, for every sort $s\in S$, the terms on the left-hand side of the pairs defining 
$\Theta^{(1)}_{s}$ are always terms in $\mathrm{CH}^{(1)}_{s}[\mathrm{Pth}_{\boldsymbol{\mathcal{A}},s}]$, while the terms on the right-hand side of such pairs are not necessarily terms in that subset of $\mathrm{T}_{\Sigma^{\boldsymbol{\mathcal{A}}}}(X)_{s}$. However, as we will show in what follows, for every sort $s\in S$, $\Theta^{(1)}_{s}$ relates pairs of terms in $\mathrm{T}_{\Sigma^{\boldsymbol{\mathcal{A}}}}(X)_{s}$ that, when mapped to $\mathrm{F}_{\Sigma^{\boldsymbol{\mathcal{A}}}}(\mathbf{Pth}_{\boldsymbol{\mathcal{A}}})_{s}$ by $\mathrm{ip}^{(1,X)@}_{s}$, give rise to pairs of paths that are related by $\mathrm{Ker}(\mathrm{CH}^{(1)})_{s}$.   

This fact already holds for the terms on the left-hand side of the pairs defining $\Theta^{(1)}_{s}$. In fact, this is a direct consequence of Proposition~\ref{PIpCH}. Therefore, we will focus our attention exclusively on the terms appearing on the right-hand side of the pairs defining $\Theta^{(1)}_{s}$.

\begin{lemma}\label{LThetaSc} 
Let $s$ be a sort in $S$ and $\mathfrak{P}$ a path in $\mathrm{Pth}_{\boldsymbol{\mathcal{A}},s}$. Then
\begin{itemize}
\item[(i)] $
\mathrm{ip}^{(1,X)@}_{s}
(\mathrm{sc}^{0\mathbf{T}_{\Sigma^{\boldsymbol{\mathcal{A}}}}(X)}_{s}
(\mathrm{CH}^{(1)}_{s}(\mathfrak{P}
)))$ 
is a  path in $\mathrm{Pth}_{\boldsymbol{\mathcal{A}},s}$;
\item[(ii)] $
\mathrm{CH}^{(1)}_{s}
(\mathrm{ip}^{(1,X)@}_{s}
(\mathrm{sc}^{0\mathbf{T}_{\Sigma^{\boldsymbol{\mathcal{A}}}}(X)}_{s}
(\mathrm{CH}^{(1)}_{s}(\mathfrak{P}
))))=
\mathrm{CH}^{(1)}_{s}
(\mathrm{sc}_{s}^{0\mathbf{Pth}_{\boldsymbol{\mathcal{A}}}}
(\mathfrak{P}
)).
$
\end{itemize}
\end{lemma}

\begin{proof}
The following chain of equalities holds
\allowdisplaybreaks
\begin{align*}
\mathrm{ip}^{(1,X)@}_{s}\left(
\mathrm{sc}^{0\mathbf{T}_{\Sigma^{\boldsymbol{\mathcal{A}}}}(X)}_{s}\left(
\mathrm{CH}^{(1)}_{s}\left(
\mathfrak{P}
\right)\right)\right)
&=
\mathrm{sc}_{s}
^{0\mathbf{F}_{\Sigma^{\boldsymbol{\mathcal{A}}}}
\left(\mathbf{Pth}_{\boldsymbol{\mathcal{A}}}\right)}\left(
\mathrm{ip}^{(1,X)@}_{s}\left(
\mathrm{CH}^{(1)}_{s}\left(
\mathfrak{P}
\right)\right)\right)
\tag{1}
\\&=
\mathrm{sc}_{s}^{0\mathbf{Pth}_{\boldsymbol{\mathcal{A}}}}\left(
\mathrm{ip}^{(1,X)@}_{s}\left(
\mathrm{CH}^{(1)}_{s}\left(
\mathfrak{P}
\right)\right)\right)
\tag{2}
\\&=
\mathrm{ip}^{(1,0)\sharp}_{s}\left(
\mathrm{sc}^{(0,1)}_{s}\left(
\mathrm{ip}^{(1,X)@}_{s}\left(
\mathrm{CH}^{(1)}_{s}\left(
\mathfrak{P}
\right)\right)\right)\right)
\tag{3}
\\&=
\mathrm{ip}^{(1,0)\sharp}_{s}\left(
\mathrm{sc}^{(0,1)}_{s}\left(
\mathfrak{P}
\right)\right)
\tag{4}
\\
&=\mathrm{sc}_{s}^{0\mathbf{Pth}_{\boldsymbol{\mathcal{A}}}}\left(
\mathfrak{P}
\right).
\tag{5}
\end{align*}

The first equality holds since $\mathrm{ip}^{(1,X)@}$ is a $\Sigma^{\boldsymbol{\mathcal{A}}}$-homomorphism in virtue of Definition~\ref{DIp}; for the second equality, we note that, by Proposition~\ref{PIpCH}, the element  $\mathrm{ip}^{(1,X)@}_{s}
(\mathrm{CH}^{(1)}_{s}
(\mathfrak{P}
))$ is a path in $[\mathfrak{P}]^{}_{s}$, thus the interpretation of the operation symbol $\mathrm{sc}^{0}_{s}$ in the $\Sigma^{\boldsymbol{\mathcal{A}}}$-algebra $\mathbf{F}_{\Sigma^{\boldsymbol{\mathcal{A}}}}(\mathbf{Pth}_{\boldsymbol{\mathcal{A}}})$ is given by the respective interpretation of the operation symbol $\mathrm{sc}^{0}_{s}$ in the $\Sigma^{\boldsymbol{\mathcal{A}}}$-algebra $\mathbf{Pth}_{\boldsymbol{\mathcal{A}}}$; the third equality  unravels the interpretation of the operation symbol $\mathrm{sc}^{0}_{s}$ in the $\Sigma^{\boldsymbol{\mathcal{A}}}$-algebra $\mathbf{Pth}_{\boldsymbol{\mathcal{A}}}$; the fourth equality follows from the fact that, by Proposition~\ref{PIpCH}, $\mathrm{ip}^{(1,X)@}_{s}
(\mathrm{CH}^{(1)}_{s}
(\mathfrak{P}
))$ is a path in $[\mathfrak{P}]^{}_{s}$, thus by Lemma~\ref{LCH} the paths $\mathrm{ip}^{(1,X)@}_{s}
(\mathrm{CH}^{(1)}_{s}
(\mathfrak{P}
))$ and $\mathfrak{P}$ have the same $(0,1)$-source; finally, the last equality unravels the interpretation of the operation symbol $\mathrm{sc}^{0}_{s}$ in the $\Sigma^{\boldsymbol{\mathcal{A}}}$-algebra $\mathbf{Pth}_{\boldsymbol{\mathcal{A}}}$.

We note that, for a sort $s\in S$ and for a path $\mathfrak{P}$ in $\mathrm{Pth}_{\boldsymbol{\mathcal{A}},s}$ it is always the case that $\mathrm{sc}_{s}^{0\mathbf{Pth}_{\boldsymbol{\mathcal{A}}}}
(\mathfrak{P})$ is a path in $\mathrm{Pth}_{\boldsymbol{\mathcal{A}},s}$.

From the just given equality,
it immediately follows that 
$$
\mathrm{CH}^{(1)}_{s}\left(
\mathrm{ip}^{(1,X)@}_{s}\left(
\mathrm{sc}^{0\mathbf{T}_{\Sigma^{\boldsymbol{\mathcal{A}}}}(X)}_{s}\left(
\mathrm{CH}^{(1)}_{s}\left(
\mathfrak{P}
\right)\right)\right)\right)
=
\mathrm{CH}^{(1)}_{s}\left(
\mathrm{sc}_{s}^{0\mathbf{Pth}_{\boldsymbol{\mathcal{A}}}}\left(
\mathfrak{P}
\right)\right).
$$
This completes the proof.
\end{proof}

A similar result holds for the $0$-target.

\begin{lemma}\label{LThetaTg}
 Let $s$ be a sort in $S$ and $\mathfrak{P}$ a path in $\mathrm{Pth}_{\boldsymbol{\mathcal{A}},s}$ Then
\begin{itemize}
\item[(i)] $
\mathrm{ip}^{(1,X)@}_{s}
(\mathrm{tg}^{0\mathbf{T}_{\Sigma^{\boldsymbol{\mathcal{A}}}}(X)}_{s}
(\mathrm{CH}^{(1)}_{s}(\mathfrak{P}
)))$ 
is a  path in $\mathrm{Pth}_{\boldsymbol{\mathcal{A}},s}$;
\item[(ii)] $
\mathrm{CH}^{(1)}_{s}
(\mathrm{ip}^{(1,X)@}_{s}
(\mathrm{tg}^{0\mathbf{T}_{\Sigma^{\boldsymbol{\mathcal{A}}}}(X)}_{s}
(\mathrm{CH}^{(1)}_{s}(\mathfrak{P}
))))=
\mathrm{CH}^{(1)}_{s}
(\mathrm{tg}_{s}^{0\mathbf{Pth}_{\boldsymbol{\mathcal{A}}}}
(\mathfrak{P}
))
$.
\end{itemize}
\end{lemma}
\begin{proof}
As it was the case in Lemma~\ref{LThetaSc}, it is easy to see that 
$$
\mathrm{ip}^{(1,X)@}_{s}\left(
\mathrm{tg}^{0\mathbf{T}_{\Sigma^{\boldsymbol{\mathcal{A}}}}(X)}_{s}\left(
\mathrm{CH}^{(1)}_{s}\left(
\mathfrak{P}
\right)\right)\right)=
\mathrm{tg}_{s}^{0\mathbf{Pth}_{\boldsymbol{\mathcal{A}}}}
\left(\mathfrak{P}
\right).
$$
The proof is analogous to that presented in Lemma~\ref{LThetaSc}. 
\end{proof}

Now we consider the case of the $0$-composition.

\begin{lemma}\label{LThetaComp} 
Let $s$ be a sort in $S$ and $\mathfrak{Q},\mathfrak{P}$ paths in $\mathrm{Pth}_{\boldsymbol{\mathcal{A}},s}$ satisfying that $\mathrm{sc}^{(0,1)}_{s}(\mathfrak{Q})=\mathrm{tg}^{(0,1)}_{s}(\mathfrak{P})$. Then
\begin{itemize}
\item[(i)] 
$\mathrm{ip}^{(1,X)@}_{s}
(\mathrm{CH}^{(1)}_{s}
(\mathfrak{Q})
\circ^{0\mathbf{T}_{\Sigma^{\boldsymbol{\mathcal{A}}}}(X)}_{s}
\mathrm{CH}^{(1)}_{s}
(\mathfrak{P}
))$ 
is a path in $\mathrm{Pth}_{\boldsymbol{\mathcal{A}},s}$;
\item[(ii)] $
\mathrm{CH}^{(1)}_{s}
(\mathrm{ip}^{(1,X)@}_{s}
(\mathrm{CH}^{(1)}_{s}
(\mathfrak{Q})
\circ^{0\mathbf{T}_{\Sigma^{\boldsymbol{\mathcal{A}}}}(X)}_{s}
\mathrm{CH}^{(1)}_{s}
(\mathfrak{P}
)))
=
\mathrm{CH}^{(1)}_{s}
(\mathfrak{Q}
\circ_{s}^{0\mathbf{Pth}_{\boldsymbol{\mathcal{A}}}}
\mathfrak{P}
)
$.
\end{itemize}
\end{lemma}

\begin{proof}
The following chain of equalities holds
\begin{flushleft}
$\mathrm{ip}^{(1,X)@}_{s}\left(
\mathrm{CH}^{(1)}_{s}\left(
\mathfrak{Q}
\right)
\circ^{0\mathbf{T}_{\Sigma^{\boldsymbol{\mathcal{A}}}}(X)}_{s}
\mathrm{CH}^{(1)}_{s}\left(
\mathfrak{P}\right)
\right)$
\allowdisplaybreaks
\begin{align*}
\quad
&=
\mathrm{ip}^{(1,X)@}_{s}\left(
\mathrm{CH}^{(1)}_{s}\left(
\mathfrak{Q}
\right)\right)
\circ^{0\mathbf{F}_{\Sigma^{\boldsymbol{\mathcal{A}}}}
(\mathbf{Pth}_{\boldsymbol{\mathcal{A}}})}_{s}
\mathrm{ip}^{(1,X)@}_{s}\left(
\mathrm{CH}^{(1)}_{s}\left(
\mathfrak{P}
\right)\right)
\tag{1}
\\&=
\mathrm{ip}^{(1,X)@}_{s}\left(
\mathrm{CH}^{(1)}_{s}\left(
\mathfrak{Q}
\right)\right)
\circ^{0\mathbf{Pth}_{\boldsymbol{\mathcal{A}}}}_{s}
\mathrm{ip}^{(1,X)@}_{s}\left(
\mathrm{CH}^{(1)}_{s}\left(
\mathfrak{P}
\right)\right)
\tag{2}
\end{align*}
\end{flushleft}

The first equality holds since  $\mathrm{ip}^{(1,X)@}$ is a $\Sigma^{\boldsymbol{\mathcal{A}}}$-homomorphism in virtue of Definition~\ref{DIp}; for the second equality, we note that, by Proposition~\ref{PIpCH}, the element $\mathrm{ip}^{(1,X)@}_{s}
(\mathrm{CH}^{(1)}_{s}
(\mathfrak{Q}
))$ is a path in $[\mathfrak{Q}]^{}_{s}$ and the element $\mathrm{ip}^{(1,X)@}_{s}
(\mathrm{CH}^{(1)}_{s}
(\mathfrak{P}
))$ is a path in $[\mathfrak{P}]^{}_{s}$. Moreover, by Lemma~\ref{LCH}, the $(0,1)$-source of $\mathrm{ip}^{(1,X)@}_{s}
(\mathrm{CH}^{(1)}_{s}
(\mathfrak{Q}
))$ is the same as the $(0,1)$-source of $\mathfrak{Q}$ and the $(0,1)$-target of $\mathrm{ip}^{(1,X)@}_{s}
(\mathrm{CH}^{(1)}_{s}
(\mathfrak{P}
))$ is the same as the $(0,1)$-target of $\mathfrak{P}$. So considering the foregoing, we have that
\allowdisplaybreaks
\begin{align*}
\mathrm{sc}^{(0,1)}_{s}\left(
\mathrm{ip}^{(1,X)@}_{s}\left(
\mathrm{CH}^{(1)}_{s}\left(
\mathfrak{Q}
\right)\right)\right)
&=
\mathrm{sc}^{(0,1)}_{s}\left(
\mathfrak{Q}
\right)
\\&=
\mathrm{tg}^{(0,1)}_{s}\left(
\mathfrak{P}
\right)
=
\mathrm{tg}^{(0,1)}_{s}\left(
\mathrm{ip}^{(1,X)@}_{s}\left(
\mathrm{CH}^{(1)}_{s}\left(
\mathfrak{P}
\right)\right)\right).
\end{align*}
Therefore, according to Proposition~\ref{PPthCatAlg}, the $0$-composition of the paths
$\mathrm{ip}^{(1,X)@}_{s}
(\mathrm{CH}^{(1)}_{s}
(\mathfrak{Q}
))$ and $\mathrm{ip}^{(1,X)@}_{s}
(\mathrm{CH}^{(1)}_{s}
(\mathfrak{P}
))$
is granted to exist in the $\Sigma^{\boldsymbol{\mathcal{A}}}$-algebra $\mathbf{Pth}_{\boldsymbol{\mathcal{A}}}$, thus the interpretation of the operation symbol $\circ^{0}_{s}$ in the $\Sigma^{\boldsymbol{\mathcal{A}}}$-algebra $\mathbf{F}_{\Sigma^{\boldsymbol{\mathcal{A}}}}(\mathbf{Pth}_{\boldsymbol{\mathcal{A}}})$ is given by the respective interpretation of the operation symbol $\circ^{0}_{s}$ in the $\Sigma^{\boldsymbol{\mathcal{A}}}$-algebra $\mathbf{Pth}_{\boldsymbol{\mathcal{A}}}$.

We note that, for the paths $\mathrm{ip}^{(1,X)@}_{s}
(\mathrm{CH}^{(1)}_{s}
(\mathfrak{Q}
))$ and $\mathrm{ip}^{(1,X)@}_{s}
(\mathrm{CH}^{(1)}_{s}
(\mathfrak{P}
))$, since its $0$-composition is granted to exist, then its $0$-composition is a path in $\mathrm{Pth}_{\boldsymbol{\mathcal{A}},s}$.

The second item in the statement of this lemma is done by Artinian induction on $(\coprod\mathrm{Pth}_{\boldsymbol{\mathcal{A}}}, \leq_{\mathbf{Pth}_{\boldsymbol{\mathcal{A}}}})$. But before considering the general case, we first consider the case in which one of the paths is a $(1,0)$-identity path.

\begin{claim}\label{CThetaComp} 
Let $s$ be a sort in $S$ and $\mathfrak{Q},\mathfrak{P}$ paths in $\mathrm{Pth}_{\boldsymbol{\mathcal{A}},s}$ such that $\mathrm{sc}^{(0,1)}_{s}(\mathfrak{Q})=\mathrm{tg}^{(0,1)}_{s}(\mathfrak{P})$. If either $\mathfrak{P}$ or $\mathfrak{Q}$ is a $(1,0)$-identity path, then
$$
\mathrm{CH}^{(1)}_{s}\left(
\mathrm{ip}^{(1,X)@}_{s}\left(
\mathrm{CH}^{(1)}_{s}\left(
\mathfrak{Q}
\right)
\circ^{0\mathbf{T}_{\Sigma^{\boldsymbol{\mathcal{A}}}}(X)}_{s}
\mathrm{CH}^{(1)}_{s}\left(
\mathfrak{P}
\right)
\right)\right)
=
\mathrm{CH}^{(1)}_{s}\left(
\mathfrak{Q}
\circ_{s}^{0\mathbf{Pth}_{\boldsymbol{\mathcal{A}}}}
\mathfrak{P}
\right)
.$$
\end{claim}

Assume without loss of generality that $\mathfrak{P}$ is a $(1,0)$-identity path. 
Note that in this case the following chain of equalities holds
\begin{flushleft}
$
\mathrm{CH}^{(1)}_{s}\left(
\mathrm{ip}^{(1,X)@}_{s}\left(
\mathrm{CH}^{(1)}_{s}\left(
\mathfrak{Q}
\right)
\circ^{0\mathbf{T}_{\Sigma^{\boldsymbol{\mathcal{A}}}}(X)}_{s}
\mathrm{CH}^{(1)}_{s}\left(
\mathfrak{P}\right)
\right)\right)$
\allowdisplaybreaks
\begin{align*}
\quad&=
\mathrm{CH}^{(1)}_{s}\left(
\mathrm{ip}^{(1,X)@}_{s}\left(
\mathrm{CH}^{(1)}_{s}\left(
\mathfrak{Q}
\right)\right)
\circ^{0\mathbf{Pth}_{\boldsymbol{\mathcal{A}}}}_{s}
\mathrm{ip}^{(1,X)@}_{s}\left(
\mathrm{CH}^{(1)}_{s}\left(
\mathfrak{P}
\right)\right)\right)
\tag{1}
\\&=
\mathrm{CH}^{(1)}_{s}\left(
\mathrm{ip}^{(1,X)@}_{s}\left(
\mathrm{CH}^{(1)}_{s}\left(
\mathfrak{Q}
\right)\right)
\circ^{0\mathbf{Pth}_{\boldsymbol{\mathcal{A}}}}_{s}
\mathfrak{P}
\right)
\tag{2}
\\&=
\mathrm{CH}^{(1)}_{s}\left(
\mathrm{ip}^{(1,X)@}_{s}\left(
\mathrm{CH}^{(1)}_{s}\left(
\mathfrak{Q}
\right)\right)\right)
\tag{3}
\\&=
\mathrm{CH}^{(1)}_{s}\left(
\mathfrak{Q}
\right)
\tag{4}
\\&=
\mathrm{CH}^{(1)}_{s}\left(
\mathfrak{Q}
\circ^{0\mathbf{Pth}_{\boldsymbol{\mathcal{A}}}}_{s}
\mathfrak{P}
\right).
\tag{5}
\end{align*}
\end{flushleft}

The first equality is a direct consequence of the equation proven in the first item of this lemma; the second equality follows from Corollary~\ref{CIpCHOneStep}, since $\mathfrak{P}$ is a $(1,0)$-identity path, we have that  $\mathrm{ip}^{(1,X)@}_{s}
(\mathrm{CH}^{(1)}_{s}
(\mathfrak{P}
))=\mathfrak{P}$; the third equality follows from the fact that $\mathfrak{P}$ is a $(1,0)$-identity path, by assumption $\mathrm{sc}^{(1,0)}_{s}(\mathfrak{Q})=\mathrm{tg}^{(1,0)}_{s}(\mathfrak{P})$, thus, by Proposition~\ref{PPthId}, $\mathfrak{P}$ is the $(1,0)$-identity path on $\mathrm{sc}^{(1,0)}_{s}(\mathfrak{Q})$. Moreover, by Proposition~\ref{PIpCH}, the element $\mathrm{ip}^{(1,X)@}_{s}
(\mathrm{CH}^{(1)}_{s}
(\mathfrak{Q}
))$ is a path in $[\mathfrak{Q}]^{}_{s}$ and, by Lemma~\ref{LCH}, the $(0,1)$-source of $\mathrm{ip}^{(1,X)@}_{s}
(\mathrm{CH}^{(1)}_{s}
(\mathfrak{Q}
))$ is the same as the $(0,1)$-source of $\mathfrak{Q}$. All in all, we have that 
$
\mathrm{ip}^{(1,X)@}_{s}
(\mathrm{CH}^{(1)}_{s}
(\mathfrak{Q}
))
\circ^{0\mathbf{Pth}_{\boldsymbol{\mathcal{A}}}}_{s}
\mathfrak{P}
=\mathrm{ip}^{(1,X)@}_{s}
(\mathrm{CH}^{(1)}_{s}
(\mathfrak{Q}
));
$ 
the fourth equality follows from the fact that, by Proposition~\ref{PIpCH}, the element $\mathrm{ip}^{(1,X)@}_{s}
(\mathrm{CH}^{(1)}_{s}
(\mathfrak{Q}
))$ is a path in $[\mathfrak{Q}]^{}_{s}$ thus its image under the Curry-Howard mapping is the same as that of $\mathfrak{Q}$; finally, the last equality follows from the fact that $\mathfrak{P}$ is the $(1,0)$-identity path on the $(0,1)$-source of $\mathfrak{Q}$ as we have mentioned before.

The same argument applies in case $\mathfrak{Q}$ is a $(1,0)$-identity path. In this case, we will argue taking into consideration that $\mathfrak{Q}$ must be the $(1,0)$-identity path on the $(0,1)$-target of $\mathfrak{P}$.

This completes the proof of Claim~\ref{CThetaComp}.

We now prove the general case by Artinian induction on $(\coprod\mathrm{Pth}_{\boldsymbol{\mathcal{A}}}, \leq_{\mathbf{Pth}_{\boldsymbol{\mathcal{A}}}})$.

\textsf{Base step of the Artinian induction}

Let $(\mathfrak{Q}\circ^{0\mathbf{Pth}_{\boldsymbol{\mathcal{A}}}}_{s}\mathfrak{P},s)$ be a minimal element in $(\coprod\mathrm{Pth}_{\boldsymbol{\mathcal{A}}},\leq_{\mathbf{Pth}_{\boldsymbol{\mathcal{A}}}})$. Then by Proposition~\ref{PMinimal}, the path $\mathfrak{Q}\circ^{0\mathbf{Pth}_{\boldsymbol{\mathcal{A}}}}_{s}\mathfrak{P}$ is either a $(1,0)$-identity path on a simple term, or an echelon.

In any case, either $\mathfrak{P}$ o $\mathfrak{Q}$ must be an $(1,0)$-identity path. The statement follows by Claim~\ref{CThetaComp}.

\textsf{Inductive step of the Artinian induction}

Let $(\mathfrak{Q}\circ^{0\mathbf{Pth}_{\boldsymbol{\mathcal{A}}}}_{s}\mathfrak{P},s)$ be a non-minimal element in $(\coprod\mathrm{Pth}_{\boldsymbol{\mathcal{A}}},\leq_{\mathbf{Pth}_{\boldsymbol{\mathcal{A}}}})$. Let us suppose that, for every sort $t\in S$ and every path $\mathfrak{Q}'\circ^{0\mathbf{Pth}_{\boldsymbol{\mathcal{A}}}}_{t}\mathfrak{P}'$ in $\mathrm{Pth}_{\boldsymbol{\mathcal{A}},t}$, if 
$
(
\mathfrak{Q}'
\circ^{0\mathbf{Pth}_{\boldsymbol{\mathcal{A}}}}_{t}
\mathfrak{P}',
t
)
<_{\mathbf{Pth}_{\boldsymbol{\mathcal{A}}}}
(
\mathfrak{Q}
\circ^{0\mathbf{Pth}_{\boldsymbol{\mathcal{A}}}}_{s}
\mathfrak{P},
s
)
$, then the statement holds for $\mathfrak{Q}'\circ^{0\mathbf{Pth}_{\boldsymbol{\mathcal{A}}}}_{t}\mathfrak{P}'$, i.e., the following equation holds
$$
\mathrm{CH}^{(1)}_{t}\left(
\mathrm{ip}^{(1,X)@}_{t}\left(
\mathrm{CH}^{(1)}_{t}\left(
\mathfrak{Q}'
\right)
\circ^{0\mathbf{T}_{\Sigma^{\boldsymbol{\mathcal{A}}}}(X)}_{t}
\mathrm{CH}^{(1)}_{t}\left(
\mathfrak{P}'
\right)\right)\right)
=
\mathrm{CH}^{(1)}_{t}\left(
\mathfrak{Q}'
\circ_{t}^{0\mathbf{Pth}_{\boldsymbol{\mathcal{A}}}}
\mathfrak{P}'
\right).
$$

Since $(\mathfrak{Q}\circ^{0\mathbf{Pth}_{\boldsymbol{\mathcal{A}}}}_{s}\mathfrak{P},s)$ is a non-minimal element in $(\coprod\mathrm{Pth}_{\boldsymbol{\mathcal{A}}},\leq_{\mathbf{Pth}_{\boldsymbol{\mathcal{A}}}})$ and taking into account, in virtue of Claim~\ref{CThetaComp}, that neither $\mathfrak{Q}$ nor $\mathfrak{P}$ are $(1,0)$-identity paths, then, by Lemma~\ref{LOrdI}, $\mathfrak{Q}\circ^{0\mathbf{Pth}_{\boldsymbol{\mathcal{A}}}}_{s}\mathfrak{P}$ is either (1) a path of length strictly greater than one containing at least one echelon or (2) an echelonless path.

If (1), then let $i\in\bb{\mathfrak{Q}\circ^{0\mathbf{Pth}_{\boldsymbol{\mathcal{A}}}}_{s}\mathfrak{P}}$ be the first index for which the one-step subpath $(\mathfrak{Q}\circ^{0\mathbf{Pth}_{\boldsymbol{\mathcal{A}}}}_{s}\mathfrak{P})^{i,i}$ of $\mathfrak{Q}\circ^{0\mathbf{Pth}_{\boldsymbol{\mathcal{A}}}}_{s}\mathfrak{P}$ is an echelon.  We distinguish the following cases according to the nature of $i$; either (1.1) $i=0$ or (1.2) $i>0$.

If~(1.1), i.e., if $\mathfrak{Q}\circ^{0\mathbf{Pth}_{\boldsymbol{\mathcal{A}}}}_{s}\mathfrak{P}$ is a path of length strictly greater than one containing an echelon on its first step, since we are assuming that $\mathfrak{P}$ is not a $(1,0)$-identity path, then $\mathfrak{P}$ has an echelon on its first step. Then it could be the case that either (1.1.1) $\mathfrak{P}$ is an echelon or (1.1.2) $\mathfrak{P}$ is a path of length strictly greater than one containing an echelon on its first step.

If~(1.1.1), i.e., if $\mathfrak{Q}\circ^{0\mathbf{Pth}_{\boldsymbol{\mathcal{A}}}}_{s}\mathfrak{P}$ is a path of length strictly greater than one containing an echelon on its first step and $\mathfrak{P}$ is an echelon then, according to Definition~\ref{DCH} the value of the Curry-Howard mapping at $\mathfrak{Q}
\circ^{0\mathbf{Pth}_{\boldsymbol{\mathcal{A}}}}_{s}
\mathfrak{P}$ is given by
$$
\mathrm{CH}^{(1)}_{s}\left(
\mathfrak{Q}
\circ^{0\mathbf{Pth}_{\boldsymbol{\mathcal{A}}}}_{s}
\mathfrak{P}
\right)
=
\mathrm{CH}^{(1)}_{s}\left(
\mathfrak{Q}\right)
\circ_{s}^{0\mathbf{T}_{\Sigma^{\boldsymbol{\mathcal{A}}}}(X)}
\mathrm{CH}^{(1)}_{s}\left(
\mathfrak{P}\right).
$$

On the other hand, the following chain of equalities holds
\begin{flushleft}
$
\mathrm{CH}^{(1)}_{s}\left(
\mathrm{ip}^{(1,X)@}_{s}\left(
\mathrm{CH}^{(1)}_{s}\left(
\mathfrak{Q}
\right)
\circ^{0\mathbf{T}_{\Sigma^{\boldsymbol{\mathcal{A}}}}(X)}_{s}
\mathrm{CH}^{(1)}_{s}\left(
\mathfrak{P}
\right)\right)\right)$
\allowdisplaybreaks
\begin{align*}
\quad
&=
\mathrm{CH}^{(1)}_{s}\left(
\mathrm{ip}^{(1,X)@}_{s}\left(
\mathrm{CH}^{(1)}_{s}\left(
\mathfrak{Q}
\right)\right)
\circ^{0\mathbf{Pth}_{\boldsymbol{\mathcal{A}}}}_{s}
\mathrm{ip}^{(1,X)@}_{s}\left(
\mathrm{CH}^{(1)}_{s}\left(
\mathfrak{P}
\right)\right)
\right)
\tag{1}
\\&=
\mathrm{CH}^{(1)}_{s}\left(
\mathrm{ip}^{(1,X)@}_{s}\left(
\mathrm{CH}^{(1)}_{s}\left(
\mathfrak{Q}
\right)\right)
\circ^{0\mathbf{Pth}_{\boldsymbol{\mathcal{A}}}}_{s}
\mathfrak{P}
\right)
\tag{2}
\\&=
\mathrm{CH}^{(1)}_{s}\left(
\mathrm{ip}^{(1,X)@}_{s}\left(
\mathrm{CH}^{(1)}_{s}\left(
\mathfrak{Q}
\right)\right)
\circ^{0\mathbf{T}_{\Sigma^{\boldsymbol{\mathcal{A}}}}(X)}_{s}
\mathrm{CH}^{(1)}_{s}\left(
\mathfrak{P}
\right)
\right)
\tag{3}
\\&=
\mathrm{CH}^{(1)}_{s}\left(
\mathfrak{Q}
\right)
\circ^{0\mathbf{T}_{\Sigma^{\boldsymbol{\mathcal{A}}}}(X)}_{s}
\mathrm{CH}^{(1)}_{s}\left(
\mathfrak{P}
\right)
\tag{4}
\\&=
\mathrm{CH}^{(1)}_{s}\left(
\mathfrak{Q}
\circ^{0\mathbf{Pth}_{\boldsymbol{\mathcal{A}}}}_{s}
\mathfrak{P}
\right).
\tag{5}
\end{align*}
\end{flushleft}

The first equality is a direct consequence of the equation proved in the first item of this lemma; the second equality follows from Corollary~\ref{CIpCHOneStep}, since by assuming that $\mathfrak{P}$ is an echelon, we have that 
$
\mathrm{ip}^{(1,X)@}_{s}
(\mathrm{CH}^{(1)}_{s}
(\mathfrak{P}
))
=\mathfrak{P}
$. Let us recall that the $0$-composition 
$\mathrm{ip}^{(1,X)@}_{s}
(\mathrm{CH}^{(1)}_{s}
(\mathfrak{Q}
))
\circ^{0\mathbf{Pth}_{\boldsymbol{\mathcal{A}}}}_{s}
\mathfrak{P}$ is well-defined since, by assumption $\mathrm{sc}^{(0,1)}_{s}(\mathfrak{Q})=\mathrm{tg}^{(0,1)}_{s}(\mathfrak{P})$, by Proposition~\ref{PIpCH} the path $\mathrm{ip}^{(1,X)@}_{s}
(\mathrm{CH}^{(1)}_{s}
(\mathfrak{Q}
))$ is in $[\mathfrak{Q}]^{}_{s}$ and, by Lemma~\ref{LCH} the $(0,1)$-source of $\mathrm{ip}^{(1,X)@}_{s}
(\mathrm{CH}^{(1)}_{s}
(\mathfrak{Q}
))$ is the same as that of $\mathfrak{Q}$; the third equality follows from the fact that  we are assuming that $\mathfrak{P}$ is an echelon, thus the path $\mathrm{ip}^{(1,X)@}_{s}
(\mathrm{CH}^{(1)}_{s}
(\mathfrak{Q}
))
\circ^{0\mathbf{Pth}_{\boldsymbol{\mathcal{A}}}}_{s}
\mathfrak{P}$ has length strictly greater than one and an echelon on its first step; the fourth equality follows from the fact that, by Proposition~\ref{PIpCH} the path $\mathrm{ip}^{(1,X)@}_{s}
(\mathrm{CH}^{(1)}_{s}
(\mathfrak{Q}
))$ is in $[\mathfrak{Q}]^{}_{s}$, thus the value of the Curry-Howard mapping at $\mathrm{ip}^{(1,X)@}_{s}
(\mathrm{CH}^{(1)}_{s}
(\mathfrak{Q}
))$ is the same as that of $\mathfrak{Q}$; finally, the last equality follows from the previous discussion on the value of the Curry-Howard mapping at $\mathfrak{Q}
\circ^{0\mathbf{Pth}_{\boldsymbol{\mathcal{A}}}}_{s}
\mathfrak{P}$.

The case of $\mathfrak{P}$ being an echelon follows.

If~(1.1.2), i.e., if $\mathfrak{Q}\circ^{0\mathbf{Pth}_{\boldsymbol{\mathcal{A}}}}_{s}\mathfrak{P}$ is a path containing an echelon on its first step and $\mathfrak{P}$ is a path of length strictly greater than one containing an echelon on its first step then, according to Definition~\ref{DCH}, the value of the Curry-Howard mapping at $\mathfrak{P}$ is given by
$$
\mathrm{CH}^{(1)}_{s}\left(\mathfrak{P}\right)
=
\mathrm{CH}^{(1)}_{s}\left(
\mathfrak{P}^{1,\bb{\mathfrak{P}}-1}
\right)
\circ_{s}^{0\mathbf{T}_{\Sigma^{\boldsymbol{\mathcal{A}}}}(X)}
\mathrm{CH}^{(1)}_{s}\left(
\mathfrak{P}^{0,0}
\right).
$$

Moreover, we have that $$\left(
\mathfrak{Q}
\circ^{0\mathbf{Pth}_{\boldsymbol{\mathcal{A}}}}_{s}
\mathfrak{P}\right)^{1,\bb{\mathfrak{Q}
\circ^{0\mathbf{Pth}_{\boldsymbol{\mathcal{A}}}}_{s}
\mathfrak{P}}-1}=\mathfrak{Q}
\circ^{0\mathbf{Pth}_{\boldsymbol{\mathcal{A}}}}_{s}
\mathfrak{P}^{1,\bb{\mathfrak{P}}-1}.$$

Thus, the value of the Curry-Howard mapping at $\mathfrak{Q}
\circ^{0\mathbf{Pth}_{\boldsymbol{\mathcal{A}}}}_{s}
\mathfrak{P}$ is given by
$$
\mathrm{CH}^{(1)}_{s}
\left(\mathfrak{Q}
\circ^{0\mathbf{Pth}_{\boldsymbol{\mathcal{A}}}}_{s}
\mathfrak{P}
\right)=
\mathrm{CH}^{(1)}_{s}\left(
\mathfrak{Q}
\circ^{0\mathbf{Pth}_{\boldsymbol{\mathcal{A}}}}_{s}
\mathfrak{P}^
{1,
\bb{
\mathfrak{P}
}-1
}
\right)
\circ^{0\mathbf{T}_{\Sigma^{\boldsymbol{\mathcal{A}}}}(X)}_{s}
\mathrm{CH}^{(1)}_{s}
\left(\mathfrak{P}^{0,0}
\right).
$$

On the other hand, the following chain of equalities holds
\begin{flushleft}
$
\mathrm{CH}^{(1)}_{s}
\left(\mathrm{ip}^{(1,X)@}_{s}
\left(\mathrm{CH}^{(1)}_{s}
(\mathfrak{Q})
\circ^{0\mathbf{T}_{\Sigma^{\boldsymbol{\mathcal{A}}}}(X)}_{s}
\mathrm{CH}^{(1)}_{s}
(\mathfrak{P}
)\right)\right)$
\allowdisplaybreaks
\begin{align*}
\quad
&=
\mathrm{CH}^{(1)}_{s}\left(
\mathrm{ip}^{(1,X)@}_{s}\left(
\mathrm{CH}^{(1)}_{s}\left(
\mathfrak{Q}
\right)\right)
\circ^{0\mathbf{Pth}_{\boldsymbol{\mathcal{A}}}}_{s}
\mathrm{ip}^{(1,X)@}_{s}\left(
\mathrm{CH}^{(1)}_{s}\left(
\mathfrak{P}
\right)\right)
\right)
\tag{1}
\\&=
\mathrm{CH}^{(1)}_{s}\Big(
\mathrm{ip}^{(1,X)@}_{s}\left(
\mathrm{CH}^{(1)}_{s}\left(
\mathfrak{Q}
\right)\right)
\circ^{0\mathbf{Pth}_{\boldsymbol{\mathcal{A}}}}_{s}
\\&\qquad\qquad\qquad\qquad
\left.
\mathrm{ip}^{(1,X)@}_{s}\left(
\mathrm{CH}^{(1)}_{s}\left(
\mathfrak{P}^{1,\bb{\mathfrak{P}}-1}\right)
\circ_{s}^{0\mathbf{T}_{\Sigma^{\boldsymbol{\mathcal{A}}}}(X)}
\mathrm{CH}^{(1)}_{s}\left(
\mathfrak{P}^{0,0}
\right)
\right)
\right)
\tag{2}
\\&=
\mathrm{CH}^{(1)}_{s}\left(
\mathrm{ip}^{(1,X)@}_{s}\left(
\mathrm{CH}^{(1)}_{s}\left(
\mathfrak{Q}
\right)\right)
\circ^{0\mathbf{Pth}_{\boldsymbol{\mathcal{A}}}}_{s}
\mathrm{ip}^{(1,X)@}_{s}\left(
\mathrm{CH}^{(1)}_{s}\left(
\mathfrak{P}^{1,\bb{\mathfrak{P}}-1}
\right)\right)
\right.
\circ^{0\mathbf{Pth}_{\boldsymbol{\mathcal{A}}}}_{s}
\\&\qquad\qquad\qquad\qquad\qquad\qquad\qquad\qquad\qquad\qquad
\mathrm{ip}^{(1,X)@}_{s}\left(
\mathrm{CH}^{(1)}_{s}\left(
\mathfrak{P}^{0,0}
\right)\right)
\Big)
\tag{3}
\\
&=
\mathrm{CH}^{(1)}_{s}\left(
\mathrm{ip}^{(1,X)@}_{s}\left(
\mathrm{CH}^{(1)}_{s}\left(
\mathfrak{Q}
\right)\right)
\circ^{0\mathbf{Pth}_{\boldsymbol{\mathcal{A}}}}_{s}
\mathrm{ip}^{(1,X)@}_{s}\left(
\mathrm{CH}^{(1)}_{s}\left(
\mathfrak{P}^{1,\bb{\mathfrak{P}}-1}
\right)\right)
\right)
\circ^{0\mathbf{T}_{\Sigma^{\boldsymbol{\mathcal{A}}}}(X)
}_{s}
\\&\qquad\qquad\qquad\qquad\qquad\qquad\qquad\qquad\qquad
\mathrm{CH}^{(1)}_{s}\left(
\mathrm{ip}^{(1,X)@}_{s}\left(
\mathrm{CH}^{(1)}_{s}\left(
\mathfrak{P}^{0,0}
\right)\right)
\right)
\tag{4}
\\
&=
\mathrm{CH}^{(1)}_{s}\left(
\mathrm{ip}^{(1,X)@}_{s}
\left(\mathrm{CH}^{(1)}_{s}\left(
\mathfrak{Q}
\right)
\circ^{0\mathbf{T}_{\Sigma^{\boldsymbol{\mathcal{A}}}}(X)}_{s}
\mathrm{CH}^{(1)}_{s}\left(
\mathfrak{P}^{1,\bb{\mathfrak{P}}-1}
\right)
\right)
\right)
\circ^{0\mathbf{T}_{\Sigma^{\boldsymbol{\mathcal{A}}}}(X)
}_{s}
\\&\qquad\qquad\qquad\qquad\qquad\qquad\qquad
\qquad\qquad\qquad\qquad\qquad\qquad
\mathrm{CH}^{(1)}_{s}\left(
\mathfrak{P}^{0,0}
\right)
\tag{5}
\\
&=
\mathrm{CH}^{(1)}_{s}\left(
\mathfrak{Q}
\circ^{0\mathbf{Pth}_{\boldsymbol{\mathcal{A}}}}_{s}
\mathfrak{P}^{1,\bb{\mathfrak{P}}-1}
\right)
\circ^{0\mathbf{T}_{\Sigma^{\boldsymbol{\mathcal{A}}}}(X)
}_{s}
\mathrm{CH}^{(1)}_{s}\left(
\mathfrak{P}^{0,0}
\right)
\tag{6}
\\&=
\mathrm{CH}^{(1)}_{s}\left(
\mathfrak{Q}\circ_{s}^{0\mathbf{Pth}_{\boldsymbol{\mathcal{A}}}}
\mathfrak{P}\right).
\tag{7}
\end{align*}
\end{flushleft}

The first equality is a direct consequence of the equation proven in the first item of this lemma; the second equality unravels the definition of the Curry-Howard mapping at $\mathfrak{P}$; the third equality follows from the fact that $\mathrm{ip}^{(1,X)@}$ is a $\Sigma^{\boldsymbol{\mathcal{A}}}$-homomorphism and from the fact that, by Proposition~\ref{PIpCH}, the elements $\mathrm{ip}^{(1,X)@}_{s}(\mathrm{CH}^{(1)}_{s}(\mathfrak{P}^{1,\bb{\mathfrak{P}}-1}))$ and $\mathrm{ip}^{(1,X)@}_{s}(\mathrm{CH}^{(1)}_{s}(\mathfrak{P}^{0,0}))$ are paths in $\mathrm{Pth}_{\boldsymbol{\mathcal{A}},s}$. We have omitted parenthesis because, by Proposition~\ref{PPthComp}, the $0$-composition of paths is associative; the fourth equality follows from the fact that $\mathfrak{P}^{0,0}$ is an echelon thus, by Corollary~\ref{CIpCHOneStep}, $\mathrm{ip}^{(1,X)@}_{s}(\mathrm{CH}^{(1)}_{s}(\mathfrak{P}^{0,0}))=\mathfrak{P}^{0,0}$. Also, by Proposition~\ref{PIpCH}, $\mathrm{ip}^{(1,X)@}_{s}(\mathrm{CH}^{(1)}_{s}(\mathfrak{P}^{1,\bb{\mathfrak{P}}-1}))$ is a path in $[\mathfrak{P}^{1,\bb{\mathfrak{P}}-1}]_{s}$, hence, by Lemma~\ref{LCH}, we have that the paths $\mathrm{ip}^{(1,X)@}_{s}(\mathrm{CH}^{(1)}_{s}(\mathfrak{P}^{1,\bb{\mathfrak{P}}-1}))$ and $\mathfrak{P}^{1,\bb{\mathfrak{P}}-1}$ have the same length. All in all, we can affirm that the path under consideration is a path of length strictly greater than one containing an echelon on its first step. All in all, the fourth equality simply unravels the Curry-Howard mapping at a path of such nature; the fifth equation follows, on one side, from the fact that, by Definition~\ref{DIp}, $\mathrm{ip}^{(1,X)@}$ is a $\Sigma^{\boldsymbol{\mathcal{A}}}$-homomorphism and, on the other side, from the fact that, according to Proposition~\ref{PIpCH}, $\mathrm{ip}^{(1,X)@}_{s}(\mathrm{CH}^{(1)}_{s}(\mathfrak{P}^{0,0}))$ is a path in $[\mathfrak{P}^{0,0}]_{s}$; the sixth equality follows from induction. Let us note that the pair $(\mathfrak{Q}\circ_{s}^{0\mathbf{Pth}_{\boldsymbol{\mathcal{A}}}}\mathfrak{P}^{1,\bb{\mathfrak{P}}-1},s)$ $\prec_{\mathbf{Pth}_{\boldsymbol{\mathcal{A}}}}$-precedes the pair $(\mathfrak{Q}\circ_{s}^{0\mathbf{Pth}_{\boldsymbol{\mathcal{A}}}}\mathfrak{P},s)$, thus we have that 
\begin{multline*}
\mathrm{CH}^{(1)}_{s}\left(
\mathrm{ip}^{(1,X)@}_{s}\left(
\mathrm{CH}^{(1)}_{s}\left(
\mathfrak{Q}
\right)
\circ_{s}^{0\mathbf{T}_{\Sigma^{\boldsymbol{\mathcal{A}}}}(X)}
\mathrm{CH}^{(1)}_{s}\left(
\mathfrak{P}^{1,\bb{\mathfrak{P}}-1}
\right)\right)\right)
\\=
\mathrm{CH}^{(1)}_{s}\left(
\mathfrak{Q}\circ_{s}^{0\mathbf{Pth}_{\boldsymbol{\mathcal{A}}}}\mathfrak{P}^{1,\bb{\mathfrak{P}}-1}
\right);
\end{multline*}
finally, the last equality recovers the value of the Curry-Howard mapping at the path $\mathfrak{Q}\circ_{s}^{0\mathbf{Pth}_{\boldsymbol{\mathcal{A}}}}\mathfrak{P}$.

The case of $\mathfrak{P}$ being a path of length strictly greater than one containing an echelon on its first step follows.

This concludes the case $i=0$.

Now, if~(1.2), i.e., if  $\mathfrak{Q}\circ^{0\mathbf{Pth}_{\boldsymbol{\mathcal{A}}}}_{s}\mathfrak{P}$ is a path of length strictly greater than one containing its first echelon at position $i\in\bb{\mathfrak{Q}\circ^{0\mathbf{Pth}_{\boldsymbol{\mathcal{A}}}}_{s}\mathfrak{P}}$ with $i>0$ then, since $\bb{\mathfrak{Q}\circ_{s}^{0\mathbf{Pth}_{\boldsymbol{\mathcal{A}}}}\mathfrak{P}}=\bb{\mathfrak{Q}}+\bb{\mathfrak{P}}$, it could be the case that either (1.2.1) $i\in\bb{\mathfrak{P}}$ or (1.2.2) $i\in [\bb{\mathfrak{P}},\bb{\mathfrak{Q}\circ_{s}^{0\mathbf{Pth}_{\boldsymbol{\mathcal{A}}}}}-1]$.

If~(1.2.1),  i.e., if  $\mathfrak{Q}\circ^{0\mathbf{Pth}_{\boldsymbol{\mathcal{A}}}}_{s}\mathfrak{P}$ is a path of length strictly greater than one containing its first echelon at position $i\in\bb{\mathfrak{P}}$ with $i>0$,  then $\mathfrak{P}$ is a path of length strictly greater than one containing an echelon on a step different from zero. Then, according to Definition~\ref{DCH}, the value of the Curry-Howard mapping at $\mathfrak{P}$ is given by 
$$
\mathrm{CH}^{(1)}_{s}\left(
\mathfrak{P}
\right)
=
\mathrm{CH}^{(1)}_{s}\left(
\mathfrak{P}^{i,\bb{\mathfrak{P}}-1}\right)
\circ_{s}^{0\mathbf{T}_{\Sigma^{\boldsymbol{\mathcal{A}}}}(X)}
\mathrm{CH}^{(1)}_{s}\left(
\mathfrak{P}^{0,i-1}
\right).
$$

Moreover,  we have that 
$$\left(
\mathfrak{Q}
\circ^{0\mathbf{Pth}_{\boldsymbol{\mathcal{A}}}}_{s}
\mathfrak{P}\right)^{i,\bb{\mathfrak{Q}
\circ^{0\mathbf{Pth}_{\boldsymbol{\mathcal{A}}}}_{s}
\mathfrak{P}}-1}=\mathfrak{Q}
\circ^{0\mathbf{Pth}_{\boldsymbol{\mathcal{A}}}}_{s}
\mathfrak{P}^{i,\bb{\mathfrak{P}}-1}.$$

Thus, the value of the Curry-Howard mapping at $\mathfrak{Q}
\circ^{0\mathbf{Pth}_{\boldsymbol{\mathcal{A}}}}_{s}
\mathfrak{P}$ is given by
$$
\mathrm{CH}^{(1)}_{s}\left(
\mathfrak{Q}
\circ^{0\mathbf{Pth}_{\boldsymbol{\mathcal{A}}}}_{s}
\mathfrak{P}
\right)=
\mathrm{CH}^{(1)}_{s}\left(
\mathfrak{Q}
\circ^{0\mathbf{Pth}_{\boldsymbol{\mathcal{A}}}}_{s}
\mathfrak{P}^
{i,
\bb{
\mathfrak{P}
}-1
}
\right)
\circ^{0\mathbf{T}_{\Sigma^{\boldsymbol{\mathcal{A}}}}(X)}_{s}
\mathrm{CH}^{(1)}_{s}\left(
\mathfrak{P}^{0,i-1}
\right).
$$

On the other hand, the following chain of equalities holds
\begin{flushleft}
$
\mathrm{CH}^{(1)}_{s}\left(
\mathrm{ip}^{(1,X)@}_{s}\left(
\mathrm{CH}^{(1)}_{s}\left(
\mathfrak{Q}
\right)
\circ^{0\mathbf{T}_{\Sigma^{\boldsymbol{\mathcal{A}}}}(X)}_{s}
\mathrm{CH}^{(1)}_{s}\left(
\mathfrak{P}
\right)\right)\right)$
\allowdisplaybreaks
\begin{align*}
\quad
&=
\mathrm{CH}^{(1)}_{s}\left(
\mathrm{ip}^{(1,X)@}_{s}\left(
\mathrm{CH}^{(1)}_{s}\left(
\mathfrak{Q}
\right)\right)
\circ^{0\mathbf{Pth}_{\boldsymbol{\mathcal{A}}}}_{s}
\mathrm{ip}^{(1,X)@}_{s}\left(
\mathrm{CH}^{(1)}_{s}\left(
\mathfrak{P}
\right)\right)
\right)
\tag{1}
\\&=
\mathrm{CH}^{(1)}_{s}
\Big(
\mathrm{ip}^{(1,X)@}_{s}\left(
\mathrm{CH}^{(1)}_{s}\left(
\mathfrak{Q}
\right)\right)
\circ^{0\mathbf{Pth}_{\boldsymbol{\mathcal{A}}}}_{s}
\\&\qquad\qquad\qquad\qquad
\mathrm{ip}^{(1,X)@}_{s}\left(
\mathrm{CH}^{(1)}_{s}\left(
\mathfrak{P}^{i,\bb{\mathfrak{P}}-1}
\right)
\circ_{s}^{0\mathbf{T}_{\Sigma^{\boldsymbol{\mathcal{A}}}}(X)}
\mathrm{CH}^{(1)}_{s}\left(
\mathfrak{P}^{0,i-1}
\right)
\right)
\Big)
\tag{2}
\\&=
\mathrm{CH}^{(1)}_{s}
\Big(
\mathrm{ip}^{(1,X)@}_{s}\left(
\mathrm{CH}^{(1)}_{s}\left(
\mathfrak{Q}
\right)\right)
\circ^{0\mathbf{Pth}_{\boldsymbol{\mathcal{A}}}}_{s}
\mathrm{ip}^{(1,X)@}_{s}\left(
\mathrm{CH}^{(1)}_{s}\left(
\mathfrak{P}^{i,\bb{\mathfrak{P}}-1}
\right)\right)
\circ^{0\mathbf{Pth}_{\boldsymbol{\mathcal{A}}}}_{s}
\\&\qquad\qquad\qquad\qquad\qquad\qquad\qquad\qquad\qquad\qquad
\mathrm{ip}^{(1,X)@}_{s}\left(
\mathrm{CH}^{(1)}_{s}\left(
\mathfrak{P}^{0,i-1}
\right)\right)
\Big)
\tag{3}
\\&=
\mathrm{CH}^{(1)}_{s}\left(
\mathrm{ip}^{(1,X)@}_{s}\left(
\mathrm{CH}^{(1)}_{s}\left(
\mathfrak{Q}
\right)\right)
\circ^{0\mathbf{Pth}_{\boldsymbol{\mathcal{A}}}}_{s}
\mathrm{ip}^{(1,X)@}_{s}\left(
\mathrm{CH}^{(1)}_{s}\left(
\mathfrak{P}^{i,\bb{\mathfrak{P}}-1}
\right)\right)
\right)
\circ^{0\mathbf{T}_{\Sigma^{\boldsymbol{\mathcal{A}}}}(X)}_{s}
\\&\qquad\qquad\qquad\qquad\qquad\qquad\qquad\qquad\,\,\,
\mathrm{CH}^{(1)}_{s}\left(
\mathrm{ip}^{(1,X)@}_{s}\left(
\mathrm{CH}^{(1)}_{s}\left(
\mathfrak{P}^{0,i-1}
\right)\right)
\right)
\tag{4}
\\&=
\mathrm{CH}^{(1)}_{s}\left(
\mathrm{ip}^{(1,X)@}_{s}\left(
\mathrm{CH}^{(1)}_{s}\left(
\mathfrak{Q}
\right)
\circ^{0\mathbf{T}_{\Sigma^{\boldsymbol{\mathcal{A}}}}(X)}_{s}
\mathrm{CH}^{(1)}_{s}\left(
\mathfrak{P}^{i,\bb{\mathfrak{P}}-1}
\right)
\right)
\right)
\circ^{0\mathbf{T}_{\Sigma^{\boldsymbol{\mathcal{A}}}}(X)}_{s}
\\&\qquad\qquad\qquad\qquad\qquad\qquad\qquad\qquad\qquad\qquad
\qquad\qquad\quad\,\,\,\,
\mathrm{CH}^{(1)}_{s}\left(
\mathfrak{P}^{0,i-1}
\right)
\tag{5}
\\&=
\mathrm{CH}^{(1)}_{s}\left(
\mathfrak{Q}
\circ^{0\mathbf{Pth}_{\boldsymbol{\mathcal{A}}}}_{s}
\mathfrak{P}^{i,\bb{\mathfrak{P}}-1}
\right)
\circ^{0\mathbf{T}_{\Sigma^{\boldsymbol{\mathcal{A}}}}(X)
}_{s}
\mathrm{CH}^{(1)}_{s}\left(
\mathfrak{P}^{0,i-1}
\right)
\tag{6}
\\&=
\mathrm{CH}^{(1)}_{s}\left(
\mathfrak{Q}\circ_{s}^{0\mathbf{Pth}_{\boldsymbol{\mathcal{A}}}}
\mathfrak{P}
\right).
\tag{7}
\end{align*}
\end{flushleft}

The first equality is a direct consequence of the equality which was proved in the first item of this lemma; the second equality unravels the definition of the Curry-Howard mapping at $\mathfrak{P}$; the third equality follows from the fact that, according to Definition~\ref{DIp}, $\mathrm{ip}^{(1,X)@}$ is a $\Sigma^{\boldsymbol{\mathcal{A}}}$-homomorphism and from the fact that, by Proposition~\ref{PIpCH}, the elements $\mathrm{ip}^{(1,X)@}_{s}(\mathrm{CH}^{(1)}_{s}(\mathfrak{P}^{i,\bb{\mathfrak{P}}-1}))$ and $\mathrm{ip}^{(1,X)@}_{s}(\mathrm{CH}^{(1)}_{s}(\mathfrak{P}^{0,i-1}))$ are paths in $\mathrm{Pth}_{\boldsymbol{\mathcal{A}},s}$. We have omitted parenthesis because, by Proposition~\ref{PPthComp}, the $0$-composition of paths is associative; the fourth equality follows from the fact that, since $\mathfrak{P}^{0,i-1}$ is an echelonless path and, by Proposition~\ref{PIpCH}, $\mathrm{ip}^{(1,X)@}_{s}(\mathrm{CH}^{(1)}_{s}(\mathfrak{P}^{0,i-1}))$ is a path in $[\mathfrak{P}^{0,i-1}]_{s}$, we have, by Lemma~\ref{LCHNEch}, that $\mathrm{ip}^{(1,X)@}_{s}(\mathrm{CH}^{(1)}_{s}(\mathfrak{P}^{0,i-1}))$ is an echelonless path. Moreover, by Lemma~\ref{LCH}, the lengths of $\mathrm{ip}^{(1,X)@}_{s}(\mathrm{CH}^{(1)}_{s}(\mathfrak{P}^{0,i-1}))$ and $\mathfrak{P}^{0,i-1}$ are equal. Also, by Proposition~\ref{PIpCH}, $\mathrm{ip}^{(1,X)@}_{s}(\mathrm{CH}^{(1)}_{s}(\mathfrak{P}^{i,\bb{\mathfrak{P}}-1}))$ is a path in $[\mathfrak{P}^{i,\bb{\mathfrak{P}}-1}]_{s}$, hence, either by Lemma~\ref{LCHUEch} or by Lemma~\ref{LCHEchInt}, $\mathrm{ip}^{(1,X)@}_{s}(\mathrm{CH}^{(1)}_{s}(\mathfrak{P}^{i,\bb{\mathfrak{P}}-1}))$ is a path containing an echelon on its first step. All things considered, we have that the path under consideration is a path of length strictly greater than one containing an echelon at position $i$. The equality simply unravels the Curry-Howard mapping at a path of such nature; the fifth equation follows, on one hand, from the fact that, by Definition~\ref{DIp}, $\mathrm{ip}^{(1,X)@}$ is a $\Sigma^{\boldsymbol{\mathcal{A}}}$-homomorphism and, on the other hand, from the fact that, according to Proposition~\ref{PIpCH}, $\mathrm{ip}^{(1,X)@}_{s}(\mathrm{CH}^{(1)}_{s}(\mathfrak{P}^{0,i-1}))$ is a path in $[\mathfrak{P}^{0,i-1}]_{s}$; the sixth equality follows from induction. Let us note that the pair $(\mathfrak{Q}\circ_{s}^{0\mathbf{Pth}_{\boldsymbol{\mathcal{A}}}}\mathfrak{P}^{i,\bb{\mathfrak{P}}-1},s)$ $\prec_{\mathbf{Pth}_{\boldsymbol{\mathcal{A}}}}$-precedes the pair $(\mathfrak{Q}\circ_{s}^{0\mathbf{Pth}_{\boldsymbol{\mathcal{A}}}}\mathfrak{P},s)$, thus we have that 
\begin{multline*}
\mathrm{CH}^{(1)}_{s}\left(
\mathrm{ip}^{(1,X)@}_{s}\left(
\mathrm{CH}^{(1)}_{s}\left(
\mathfrak{Q}
\right)
\circ_{s}^{0\mathbf{T}_{\Sigma^{\boldsymbol{\mathcal{A}}}}(X)}
\mathrm{CH}^{(1)}_{s}\left(
\mathfrak{P}^{i,\bb{\mathfrak{P}}-1}
\right)\right)\right)
\\=
\mathrm{CH}^{(1)}_{s}\left(
\mathfrak{Q}\circ_{s}^{0\mathbf{Pth}_{\boldsymbol{\mathcal{A}}}}\mathfrak{P}^{i,\bb{\mathfrak{P}}-1}
\right);
\end{multline*}
finally, the last equality recovers the value of the Curry-Howard mapping at  path $\mathfrak{Q}\circ_{s}^{0\mathbf{Pth}_{\boldsymbol{\mathcal{A}}}}\mathfrak{P}$.

The case $i\in\bb{\mathfrak{P}}$ follows.

If~(1.2.2), i.e., if $\mathfrak{Q}\circ^{0\mathbf{Pth}_{\boldsymbol{\mathcal{A}}}}_{s}\mathfrak{P}$ is a path of length strictly greater than one containing its first echelon at position $i\in[\bb{\mathfrak{P}},\bb{\mathfrak{Q}\circ_{s}^{0\mathbf{Pth}_{\boldsymbol{\mathcal{A}}}}\mathfrak{P}}-1]$ then $\mathfrak{Q}$ is a non-$(1,0)$-identity path containing an echelon, whilst $\mathfrak{P}$ is an echelonless path.

We will distinguish three cases according to whether (1.2.2.1) $\mathfrak{Q}$ is an echelon; (1.2.2.2) $\mathfrak{Q}$ is a path of length strictly greater than one containing an echelon on its first step or (1.2.2.3) $\mathfrak{Q}$ is a path of length strictly greater than one containing an echelon on a step different from zero. These cases can be proven following a similar argument to those three cases presented above. We leave the details for the interested reader.

If (2), i.e., if $\mathfrak{Q}
\circ^{0\mathbf{Pth}_{\boldsymbol{\mathcal{A}}}}_{s}
\mathfrak{P}$ is an echelonless path then, regarding the paths $\mathfrak{Q}$ and $\mathfrak{P}$ we have that
\begin{itemize}
\item[(i)] $\mathfrak{P}$ is an echelonless path;
\item[(ii)] $\mathfrak{Q}$ is an echelonless path.
\end{itemize}

In this case, by Lemma~\ref{LPthHeadCt}, there exists a unique word $\mathbf{s}\in S^{\star}-\{\lambda\}$ and a unique operation symbol $\sigma\in\Sigma_{\mathbf{s},s}$  such that 
the family of translations occurring in $\mathfrak{P}$ is a family of translations of type $\sigma$.

Since $\mathfrak{Q}
\circ^{0\mathbf{Pth}_{\boldsymbol{\mathcal{A}}}}_{s}
\mathfrak{P}$ is associated to the operation symbol $\sigma$, in virtue of Lemma~\ref{LPthHeadCt}, all translations appearing in $\mathfrak{Q}$ and $\mathfrak{P}$ have the same type, i.e., the paths $\mathfrak{Q}$ and $\mathfrak{P}$ are also associated to this same operation symbol $\sigma$.  

Let $((\mathfrak{Q}
\circ^{0\mathbf{Pth}_{\boldsymbol{\mathcal{A}}}}_{s}
\mathfrak{P})_{j})_{j\in\bb{\mathbf{s}}}$ be the family of paths we can extract from 
$\mathfrak{Q}
\circ^{0\mathbf{Pth}_{\boldsymbol{\mathcal{A}}}}_{s}
\mathfrak{P}$ in virtue of Lemma~\ref{LPthExtract}. Then, the value of the Curry-Howard mapping at  
$\mathfrak{Q}
\circ^{0\mathbf{Pth}_{\boldsymbol{\mathcal{A}}}}_{s}
\mathfrak{P}$ is given by
$$
\mathrm{CH}^{(1)}_{s}
\left(\mathfrak{Q}
\circ^{0\mathbf{Pth}_{\boldsymbol{\mathcal{A}}}}_{s}
\mathfrak{P}
\right)
=
\sigma
^{\mathbf{T}_{\Sigma^{\boldsymbol{\mathcal{A}}}}(X)}
\left(\left(\mathrm{CH}^{(1)}_{s_{j}}\left(\left(
\mathfrak{Q}
\circ^{0\mathbf{Pth}_{\boldsymbol{\mathcal{A}}}}_{s}
\mathfrak{P}
\right)_{j}
\right)\right)_{j\in\bb{\mathbf{s}}}
\right).
$$

The following equality is a direct consequence of the equality proved in the first item of this lemma
\begin{multline*}
\mathrm{CH}^{(1)}_{s}\left(
\mathrm{ip}^{(1,X)@}_{s}\left(
\mathrm{CH}^{(1)}_{s}\left(
\mathfrak{Q}
\right)
\circ^{0\mathbf{T}_{\Sigma^{\boldsymbol{\mathcal{A}}}}(X)}_{s}
\mathrm{CH}^{(1)}_{s}\left(
\mathfrak{P}
\right)
\right)
\right)
\\=
\mathrm{CH}^{(1)}_{s}\left(
\mathrm{ip}^{(1,X)@}_{s}\left(
\mathrm{CH}^{(1)}_{s}\left(
\mathfrak{Q}
\right)
\right)
\circ^{0\mathbf{Pth}_{\boldsymbol{\mathcal{A}}}}_{s}
\mathrm{ip}^{(1,X)@}_{s}\left(
\mathrm{CH}^{(1)}_{s}\left(
\mathfrak{P}
\right)
\right)
\right)
\end{multline*}

We will study the  nature of the path $\mathfrak{R}$ given by
$$
\mathfrak{R}=
\mathrm{ip}^{(1,X)@}_{s}\left(
\mathrm{CH}^{(1)}_{s}\left(
\mathfrak{Q}
\right)
\right)
\circ^{0\mathbf{Pth}_{\boldsymbol{\mathcal{A}}}}_{s}
\mathrm{ip}^{(1,X)@}_{s}\left(
\mathrm{CH}^{(1)}_{s}\left(
\mathfrak{P}
\right)
\right).
$$

Note that
\begin{itemize}
\item[(i)] $\mathfrak{P}$ is an echelonless path associated to the operation symbol $\sigma$.
\end{itemize}

Let $(\mathfrak{P}_{j})_{j\in\bb{\mathbf{s}}}$ be the family of paths we can extract from $\mathfrak{P}$ in virtue of Lemma~\ref{LPthExtract}. Then, according to Definition~\ref{DCH}, the value of the Curry-Howard mapping at $\mathfrak{P}$ is given by
$$
\mathrm{CH}^{(1)}_{s}\left(
\mathfrak{P}
\right)
=
\sigma^{\mathbf{T}_{\Sigma^{\boldsymbol{\mathcal{A}}}}(X)}
\left(\left(
\mathrm{CH}^{(1)}_{s_{j}}\left(
\mathfrak{P}_{j}\right)
\right)_{j\in\bb{\mathbf{s}}}
\right).
$$

By Lemma~\ref{LCHNEch}, $\mathrm{CH}^{(1)}_{s}(\mathfrak{P})\in\mathcal{T}(\sigma,\mathrm{T}_{\Sigma^{\boldsymbol{\mathcal{A}}}}(X))_{1}$, which is a subset of $\mathrm{T}_{\Sigma^{\boldsymbol{\mathcal{A}}}}(X)^{\mathsf{E}}_{s}$. By Proposition~\ref{PIpCH} $\mathrm{ip}^{(1,X)@}_{s}(\mathrm{CH}^{(1)}_{s}(\mathfrak{P}))$ is a path in $[\mathfrak{P}]_{s}$ then, by Lemma~\ref{LCHNEch}, we have that
\begin{itemize}
\item[(i)] $\mathrm{ip}^{(1,X)@}_{s}(\mathrm{CH}^{(1)}_{s}(\mathfrak{P}))$ is an echelonless path associated to $\sigma$.
\end{itemize}

Let $(\mathfrak{P}'_{j})_{j\in\bb{\mathbf{s}}}$ be the family of paths we can extract from $\mathrm{ip}^{(1,X)@}_{s}(\mathrm{CH}^{(1)}_{s}(\mathfrak{P}))$ in virtue of Lemma~\ref{LPthExtract}. Then, according to Definition~\ref{DCH}, the value of the Curry-Howard mapping at $\mathrm{ip}^{(1,X)@}_{s}(\mathrm{CH}^{(1)}_{s}(\mathfrak{P}))$ is given by
$$
\mathrm{CH}^{(1)}_{s}\left(
\mathrm{ip}^{(1,X)@}_{s}\left(
\mathrm{CH}^{(1)}_{s}\left(
\mathfrak{P}
\right)\right)
\right)
=
\sigma^{\mathbf{T}_{\Sigma^{\boldsymbol{\mathcal{A}}}}(X)}
\left(\left(
\mathrm{CH}^{(1)}_{s_{j}}\left(
\mathfrak{P}'_{j}\right)
\right)_{j\in\bb{\mathbf{s}}}
\right).
$$

By Proposition~\ref{PIpCH} $\mathrm{ip}^{(1,X)@}_{s}(\mathrm{CH}^{(1)}_{s}(\mathfrak{P}))$ is a path in $[\mathfrak{P}]_{s}$. Therefore, we have that, for every $j\in\bb{\mathbf{s}}$, it is the case that
$$
\left(
\mathfrak{P}_{j},
\mathfrak{P}'_{j}
\right)
\in\mathrm{Ker}\left(
\mathrm{CH}^{(1)}
\right)_{s_{j}}.
$$

On the other hand
\begin{itemize}
\item[(i)] $\mathfrak{Q}$ is an echelonless path associated to the operation symbol $\sigma$.
\end{itemize}

Let $(\mathfrak{Q}_{j})_{j\in\bb{\mathbf{s}}}$ be the family of paths we can extract from $\mathfrak{Q}$ in virtue of Lemma~\ref{LPthExtract}. Then, according to Definition~\ref{DCH}, the value of the Curry-Howard mapping at $\mathfrak{Q}$ is given by
$$
\mathrm{CH}^{(1)}_{s}\left(
\mathfrak{Q}
\right)
=
\sigma^{\mathbf{T}_{\Sigma^{\boldsymbol{\mathcal{A}}}}(X)}
\left(\left(
\mathrm{CH}^{(1)}_{s_{j}}\left(
\mathfrak{Q}_{j}\right)
\right)_{j\in\bb{\mathbf{s}}}
\right).
$$

By Lemma~\ref{LCHNEch}, $\mathrm{CH}^{(1)}_{s}(\mathfrak{Q})\in\mathcal{T}(\sigma,\mathrm{T}_{\Sigma^{\boldsymbol{\mathcal{A}}}}(X))_{1}$, which is a subset of $\mathrm{T}_{\Sigma^{\boldsymbol{\mathcal{A}}}}(X)^{\mathsf{E}}_{s}$. By Proposition~\ref{PIpCH} $\mathrm{ip}^{(1,X)@}_{s}(\mathrm{CH}^{(1)}_{s}(\mathfrak{Q}))$ is a path in $[\mathfrak{Q}]_{s}$ then, by Lemma~\ref{LCHNEch}, we have that
\begin{itemize}
\item[(i)] $\mathrm{ip}^{(1,X)@}_{s}(\mathrm{CH}^{(1)}_{s}(\mathfrak{Q}))$ is an echelonless path associated to $\sigma$.
\end{itemize}

Let $(\mathfrak{Q}'_{j})_{j\in\bb{\mathbf{s}}}$ be the family of paths we can extract from $\mathrm{ip}^{(1,X)@}_{s}(\mathrm{CH}^{(1)}_{s}(\mathfrak{Q}))$ in virtue of Lemma~\ref{LPthExtract}. Then, according to Definition~\ref{DCH}, the value of the Curry-Howard mapping at $\mathrm{ip}^{(1,X)@}_{s}(\mathrm{CH}^{(1)}_{s}(\mathfrak{Q}))$ is given by
$$
\mathrm{CH}^{(1)}_{s}\left(
\mathrm{ip}^{(1,X)@}_{s}\left(
\mathrm{CH}^{(1)}_{s}\left(
\mathfrak{Q}
\right)\right)
\right)
=
\sigma^{\mathbf{T}_{\Sigma^{\boldsymbol{\mathcal{A}}}}(X)}
\left(\left(
\mathrm{CH}^{(1)}_{s_{j}}\left(
\mathfrak{Q}'_{j}\right)
\right)_{j\in\bb{\mathbf{s}}}
\right).
$$

By Proposition~\ref{PIpCH} $\mathrm{ip}^{(1,X)@}_{s}(\mathrm{CH}^{(1)}_{s}(\mathfrak{Q}))$ is a path in $[\mathfrak{Q}]_{s}$. Therefore, we have that, for every $j\in\bb{\mathbf{s}}$, it is the case that
$$
\left(
\mathfrak{Q}_{j},
\mathfrak{Q}'_{j}
\right)
\in\mathrm{Ker}\left(
\mathrm{CH}^{(1)}
\right)_{s_{j}}.
$$

Taking into account items (i) and (ii), we conclude that the path $\mathfrak{R}$ given by
$$
\mathfrak{R}=
\mathrm{ip}^{(1,X)@}_{s}\left(
\mathrm{CH}^{(1)}_{s}\left(
\mathfrak{Q}
\right)
\right)
\circ^{0\mathbf{Pth}_{\boldsymbol{\mathcal{A}}}}_{s}
\mathrm{ip}^{(1,X)@}_{s}\left(
\mathrm{CH}^{(1)}_{s}\left(
\mathfrak{P}
\right)
\right).
$$
is an echelonless path associated to the operation symbol symbol $\sigma$.

Let $(\mathfrak{R}_{j})_{j\in\bb{\mathbf{s}}}$ be the family of paths we can extract from $\mathfrak{R}$ in virtue of Lemma~\ref{LPthExtract}. For every $j\in\bb{\mathbf{s}}$, the following chain of equalities holds
\begin{align*}
\mathfrak{R}_{j}&=\mathfrak{Q}'_{j}\circ^{0\mathbf{Pth}_{\boldsymbol{\mathcal{A}}}}_{s_{j}}
\mathfrak{P}'_{j}
\tag{1}
\\&=
\mathrm{ip}^{(1,X)@}_{s_{j}}\left(
\mathrm{CH}^{(1)}_{s_{j}}\left(
\mathfrak{Q}_{j}
\right)
\right)
\circ_{s_{j}}^{0\mathbf{Pth}_{\boldsymbol{\mathcal{A}}}}
\mathrm{ip}^{(1,X)@}_{s_{j}}\left(
\mathrm{CH}^{(1)}_{s_{j}}\left(
\mathfrak{P}_{j}
\right)
\right).
\tag{2}
\end{align*}

the first equality follows from the proof of Lemma~\ref{LPthExtract} and the description of $\mathfrak{R}$; finally, the last equality follows from the description of $\mathfrak{R}$ as a $0$-composition of normalized paths. Thus, the $j$-th component we can extract from the normalized path $\mathrm{ip}^{(1,X)@}_{s}(\mathrm{CH}^{(1)}_{s}(\mathfrak{Q}))$ is the normalization of the $j$-th component we can extract from $\mathfrak{Q}$. Analogously for $\mathrm{ip}^{(1,X)@}_{s}(\mathrm{CH}^{(1)}_{s}(\mathfrak{P}))$.

Also from the proof of Lemma~\ref{LPthExtract}, we also have that, for every $j\in\bb{\mathbf{s}}$,
$$
\left(
\mathfrak{Q}
\circ^{0\mathbf{Pth}_{\boldsymbol{\mathcal{A}}}}_{s}
\mathfrak{P}
\right)_{j}
=
\mathfrak{Q}_{j}
\circ^{0\mathbf{Pth}_{\boldsymbol{\mathcal{A}}}}_{s_{j}}
\mathfrak{P}_{j}
.
$$

Therefore, we are in position to conclude the initial discussion at the beginning of this subcase. The following chain of equalities holds
\begin{flushleft}
$
\mathrm{CH}^{(1)}_{s}\left(
\mathrm{ip}^{(1,X)@}_{s}\left(
\mathrm{CH}^{(1)}_{s}\left(
\mathfrak{Q}
\right)
\circ^{0\mathbf{T}_{\Sigma^{\boldsymbol{\mathcal{A}}}}(X)}_{s}
\mathrm{CH}^{(1)}_{s}\left(
\mathfrak{P}
\right)
\right)
\right)
$
\begin{align*}
\quad&=
\mathrm{CH}^{(1)}_{s}\left(
\mathrm{ip}^{(1,X)@}_{s}\left(
\mathrm{CH}^{(1)}_{s}\left(
\mathfrak{Q}
\right)
\right)
\circ^{0\mathbf{Pth}_{\boldsymbol{\mathcal{A}}}}_{s}
\mathrm{ip}^{(1,X)@}_{s}\left(
\mathrm{CH}^{(1)}_{s}\left(
\mathfrak{P}
\right)
\right)
\right)
\tag{1}
\\&=
\sigma^{\mathbf{T}_{\Sigma^{\boldsymbol{\mathcal{A}}}}(X)}
\bigg(\left(
\mathrm{CH}^{(1)}_{s_{j}}\left(
\mathrm{ip}^{(1,X)@}_{s_{j}}\left(
\mathrm{CH}^{(1)}_{s_{j}}\left(
\mathfrak{Q}_{j}
\right)
\right)
\circ_{s_{j}}^{0\mathbf{Pth}_{\boldsymbol{\mathcal{A}}}}
\right.
\right.
\\&\qquad\qquad\qquad\qquad\qquad\qquad\qquad\qquad\qquad
\mathrm{ip}^{(1,X)@}_{s_{j}}\left(
\mathrm{CH}^{(1)}_{s_{j}}\left(
\mathfrak{P}_{j}
\right)
\right)
\Big)
\Big)_{j\in\bb{\mathbf{s}}}
\bigg)
\tag{2}
\\&=
\sigma^{\mathbf{T}_{\Sigma^{\boldsymbol{\mathcal{A}}}}(X)}
\bigg(\left(
\mathrm{CH}^{(1)}_{s_{j}}\left(
\mathrm{ip}^{(1,X)@}_{s_{j}}\left(
\mathrm{CH}^{(1)}_{s_{j}}\left(
\mathfrak{Q}_{j}
\right)
\circ_{s_{j}}^{0
\mathbf{T}_{\Sigma^{\boldsymbol{\mathcal{A}}}}(X)
}
\right.
\right.
\right.
\\&\qquad\qquad\qquad\qquad\qquad\qquad\qquad\qquad\qquad\qquad\quad\,\,\,
\mathrm{CH}^{(1)}_{s_{j}}\left(
\mathfrak{P}_{j}
\right)
\Big)
\Big)
\Big)_{j\in\bb{\mathbf{s}}}
\bigg)
\tag{3}
\\&=
\sigma^{\mathbf{T}_{\Sigma^{\boldsymbol{\mathcal{A}}}}(X)}
\left(\left(
\mathrm{CH}^{(1)}_{s_{j}}\left(
\mathfrak{Q}_{j}
\circ^{0\mathbf{Pth}_{\boldsymbol{\mathcal{A}}}}_{s_{j}}
\mathfrak{P}_{j}
\right)
\right)_{j\in\bb{\mathbf{s}}}
\right)
\tag{4}
\\&=
\mathrm{CH}^{(1)}_{s}\left(
\mathfrak{Q}
\circ^{0\mathbf{Pth}_{\boldsymbol{\mathcal{A}}}}_{s}
\mathfrak{P}
\right).
\tag{5}
\end{align*}
\end{flushleft}

The first equality was already proven at the beginning of this case; the second equality simply unravels the Curry-Howard mapping at $\mathfrak{R}$; the third equality follows from the fact that, by Definition~\ref{DIp}, $\mathrm{ip}^{(1,X)@}$ is a $\Sigma^{\boldsymbol{\mathcal{A}}}$-homomorphism; the fourth equality follows from induction. Let us note that the pair $(\mathfrak{Q}_{j}\circ^{0\mathbf{Pth}_{\boldsymbol{\mathcal{A}}}}_{s_{j}}\mathfrak{P}_{j}, s_{j})$ $\prec_{\mathbf{Pth}_{\boldsymbol{\mathcal{A}}}}$-precedes $(\mathfrak{Q}\circ^{0\mathbf{Pth}_{\boldsymbol{\mathcal{A}}}}_{s}\mathfrak{P}, s)$, thus we have that
$$
\mathrm{CH}^{(1)}_{s_{j}}\left(
\mathrm{ip}^{(1,X)@}_{s_{j}}\left(
\mathrm{CH}^{(1)}_{s_{j}}\left(
\mathfrak{Q}_{j}
\right)
\circ^{0\mathbf{T}_{\Sigma^{\boldsymbol{\mathcal{A}}}}(X)}_{s_{j}}
\mathrm{CH}^{(1)}_{s_{j}}\left(
\mathfrak{P}_{j}
\right)
\right)
\right)
=
\mathrm{CH}^{(1)}_{s_{j}}\left(
\mathfrak{Q}_{j}
\circ^{0\mathbf{Pth}_{\boldsymbol{\mathcal{A}}}}_{s_{j}}
\mathfrak{P}_{j}
\right);
$$
finally, the last equality recovers the value of the Curry-Howard mapping at the path $\mathfrak{Q}
\circ^{0\mathbf{Pth}_{\boldsymbol{\mathcal{A}}}}_{s}
\mathfrak{P}$.

Case~(2) follows.

This completes the proof.
\end{proof}

In the following remark we warn the reader about the fact that $\mathrm{ip}^{(1,X)@}\circ \mathrm{CH}^{(1)}$ does not commutes, in general, with the $0$-composition.

\begin{remark}\label{RIpCHComp} 
For the set of sorts $S=1$, let $\Sigma$ be the signature that only contains a binary operation $\diamond$,  for $X=\{x,y,z\}$, let $\mathbf{T}_{\Sigma}(X)$ be the free $\Sigma$-algebra on $X$ and let $\mathcal{A}$  be the one-sorted set of rewrite rules consisting of
$$
\begin{array}{rclcrcl}
\mathfrak{p}&=&(y,z),&
\qquad&
\mathfrak{q}&=&(x,y)
\end{array}
$$

For the paths
$$
\xymatrix@C=55pt{
\mathfrak{P}: x\diamond y
\ar[r]^-{(\mathfrak{p}, x\diamond\underline{\quad})}
&
x\diamond z
},
\qquad\qquad
\xymatrix@C=55pt{
\mathfrak{Q}: x\diamond z
\ar[r]^-{(\mathfrak{q}, \underline{\quad}\diamond z)}
&
y\diamond z
},
$$
we have that  
$\mathrm{CH}^{(1)}(\mathfrak{P})=x\diamond \mathfrak{p}$ and 
$\mathrm{CH}^{(1)}(\mathfrak{Q})=\mathfrak{q}\diamond z$.
Moreover, since both $\mathfrak{P}$ and $\mathfrak{Q}$ are one-step paths, we have, by Corollary~\ref{CIpCHOneStep}, that 
$\mathrm{ip}^{(1,X)@}(\mathrm{CH}^{(1)}(\mathfrak{P}))=\mathfrak{P}$ and $\mathrm{ip}^{(1,X)@}(\mathrm{CH}^{(1)}(\mathfrak{Q}))=\mathfrak{Q}$.

Since, in addition, $\mathrm{sc}^{(0,1)}(\mathfrak{Q})=\mathrm{tg}^{(0,1)}(\mathfrak{P})$, we can consider the $0$-composition
$$
\xymatrix@C=55pt{
\mathfrak{Q}\circ^{0\mathbf{Pth}_{\boldsymbol{\mathcal{A}}}}\mathfrak{P}: x\diamond y
\ar[r]^-{(\mathfrak{p}, x\diamond\underline{\quad})}
&
x\diamond z
\ar[r]^-{(\mathfrak{q}, \underline{\quad}\diamond z)}
&
y\diamond z
},
$$

In this case, $\mathrm{CH}^{(1)}(\mathfrak{Q}\circ^{0\mathbf{Pth}_{\boldsymbol{\mathcal{A}}}}\mathfrak{P})=\mathfrak{q}\diamond \mathfrak{p}$. Moreover,
$$
\xymatrix@C=55pt{
\mathrm{ip}^{(1,X)@}(\mathrm{CH}^{(1)}(\mathfrak{Q}\circ^{0\mathbf{Pth}_{\boldsymbol{\mathcal{A}}}}\mathfrak{P})): x\diamond y
\ar[r]^-{(\mathfrak{q}, \underline{\quad}\diamond y)}
&
y\diamond y
\ar[r]^-{(\mathfrak{p}, y\diamond\underline{\quad})}
&
y\diamond z
}.
$$
Hence
$$
\mathrm{ip}^{(1,X)@}_{s}\left(
\mathrm{CH}^{(1)}_{s}\left(
\mathfrak{Q}
\circ_{s}^{0\mathbf{Pth}_{\boldsymbol{\mathcal{A}}}}
\mathfrak{P}
\right)
\right)
\neq
\mathrm{ip}^{(1,X)@}_{s}\left(
\mathrm{CH}^{(1)}_{s}\left(
\mathfrak{Q}
\right)
\right)
\circ_{s}^{0\mathbf{Pth}_{\boldsymbol{\mathcal{A}}}}
\mathrm{ip}^{(1,X)@}_{s}\left(
\mathrm{CH}^{(1)}_{s}\left(
\mathfrak{P}
\right)
\right).
$$

\end{remark}

\section{
\texorpdfstring
{The congruence $\Theta^{[1]}$ on $\mathbf{T}_{\Sigma^{\boldsymbol{\mathcal{A}}}}(X)$}
{A congruence for the free algebra}
}

We next show that the properties that hold for the relation $\Theta^{(1)}$ also hold for the smallest 
$\Sigma^{\boldsymbol{\mathcal{A}}}$-congru\-ence on $\mathbf{T}_{\Sigma^{\boldsymbol{\mathcal{A}}}}(X)$ containing $\Theta^{(1)}$. In what follows we will use the notation introduced in Definition~\ref{DCongOpInt}.

\begin{restatable}{definition}{DThetaCong}
\label{DThetaCong}
\index{Theta!first-order!$\Theta^{[1]}$}
We let $\Theta^{[1]}$ stand for $\mathrm{Cg}_{\mathbf{T}_{\Sigma^{\boldsymbol{\mathcal{A}}}}(X)}
\left(\Theta^{(1)}\right)$, the smallest $\Sigma^{\boldsymbol{\mathcal{A}}}$-congru\-ence on $\mathbf{T}_{\Sigma^{\boldsymbol{\mathcal{A}}}}(X)$ containing $\Theta^{(1)}$.
\end{restatable}

\begin{remark}\label{RThetaCong}
To simplify the notation, we let $\mathrm{C}_{\Sigma^{\boldsymbol{\mathcal{A}}}}$ stand for the congruence generating operator $\mathrm{C}_{\mathbf{T}_{\Sigma^{\boldsymbol{\mathcal{A}}}}(X)}$ on $\mathrm{T}_{\Sigma^{\boldsymbol{\mathcal{A}}}}(X)\times \mathrm{T}_{\Sigma^{\boldsymbol{\mathcal{A}}}}(X)$ determined by 
$\mathbf{T}_{\Sigma^{\boldsymbol{\mathcal{A}}}}(X)$ (see Definition~\ref{DCongOpC} for the general case).
Let us also recall that, by Definition~\ref{DCongOpC}, if $\Phi\subseteq \mathrm{T}_{\Sigma^{\boldsymbol{\mathcal{A}}}}(X)\times \mathrm{T}_{\Sigma^{\boldsymbol{\mathcal{A}}}}(X)$, then  
$$
\mathrm{C}_{\Sigma^{\boldsymbol{\mathcal{A}}}}\left(
\Phi
\right)
=\left(
\Phi\circ\Phi
\right)
\cup
\left(
\bigcup_{\tau\in\Sigma^{\boldsymbol{\mathcal{A}}}_{\neq\lambda,s}}
\tau^{\mathbf{T}_{\Sigma^{\boldsymbol{\mathcal{A}}}}(X)}
\times
\tau^{\mathbf{T}_{\Sigma^{\boldsymbol{\mathcal{A}}}}(X)}
\left[
\Phi_{\mathrm{ar}(\tau)}
\right]
\right)_{s\in S}.
$$
Moreover, for the family $(\mathrm{C}^{n}_{\Sigma^{\boldsymbol{\mathcal{A}}}}(\Theta^{(1)}))_{n\in\mathbb{N}}$ in $\mathrm{Sub}(\mathrm{T}_{\Sigma^{\boldsymbol{\mathcal{A}}}}(X)^{2})$, defined recursively, as:
\allowdisplaybreaks
\begin{align*}
\mathrm{C}^{0}_{\Sigma^{\boldsymbol{\mathcal{A}}}}\left(
\Theta^{(1)}\right)
&=
\Theta^{(1)}\cup\left(\Theta^{(1)}\right)^{-1}\cup\Delta_{\mathrm{T}_{\Sigma^{\boldsymbol{\mathcal{A}}}}(X)},\\
\mathrm{C}^{n+1}_{\Sigma^{\boldsymbol{\mathcal{A}}}}\left(
\Theta^{(1)}\right)
&=\mathrm{C}_{\Sigma^{\boldsymbol{\mathcal{A}}}}
\left(
\mathrm{C}^{n}_{\Sigma^{\boldsymbol{\mathcal{A}}}}\left(
\Theta^{(1)}\right)\right),\, n\geq 0,
\end{align*}
we have, by Proposition~\ref{PCongIntC}, that 
$$\mathrm{C}^{\omega}_{\Sigma^{\boldsymbol{\mathcal{A}}}}\left(
\Theta^{(1)}\right)
=\bigcup_{n\in\mathbb{N}}\mathrm{C}^{n}_{\Sigma^{\boldsymbol{\mathcal{A}}}}\left(
\Theta^{(1)}
\right)
=\Theta^{[1]}.
$$ 
\end{remark}

\begin{restatable}{proposition}{PThetaCongCatAlg}
\label{PThetaCongCatAlg}
\index{terms!first-order!$\mathrm{T}_{\Sigma^{\boldsymbol{\mathcal{A}}}}(X)/{\Theta^{[1]}}$}
\index{terms!first-order!$[P]_{\Theta^{[1]}_{s}}$}
The $S$-sorted set $\mathrm{T}_{\Sigma^{\boldsymbol{\mathcal{A}}}}(X)/{\Theta^{[1]}}$ is equipped, in a natural way, with a structure of $\Sigma^{\boldsymbol{\mathcal{A}}}$-algebra, we denote by 
$\mathbf{T}_{\Sigma^{\boldsymbol{\mathcal{A}}}}(X)/{\Theta^{[1]}}$ the corresponding $\Sigma^{\boldsymbol{\mathcal{A}}}$-algebra. Moreover, 
\index{projection!first-order!$\mathrm{pr}^{\Theta^{[1]}}$}
the canonical projection $\mathrm{pr}^{\Theta^{[1]}}$ from $\mathrm{T}_{\Sigma^{\boldsymbol{\mathcal{A}}}}(X)$ to $\mathrm{T}_{\Sigma^{\boldsymbol{\mathcal{A}}}}(X)/{\Theta^{[1]}}$ that, for every sort $s\in S$, sends  a term $P$ in $\mathrm{T}_{\Sigma^{\boldsymbol{\mathcal{A}}}}(X)_{s}$ to $[P]_{\Theta^{[1]}_{s}}$ in 
$\mathrm{T}_{\Sigma^{\boldsymbol{\mathcal{A}}}}(X)_{s}/{\Theta^{[1]}_{s}}$
determines a surjective $\Sigma^{\boldsymbol{\mathcal{A}}}$-homomorphism from $\mathbf{T}_{\Sigma^{\boldsymbol{\mathcal{A}}}}(X)$ to $\mathbf{T}_{\Sigma^{\boldsymbol{\mathcal{A}}}}(X)/{\Theta^{[1]}}$.
\end{restatable}

We next introduce the $\Sigma$-reduct of the $\Sigma^{\boldsymbol{\mathcal{A}}}$-algebra $\mathbf{T}_{\Sigma^{\boldsymbol{\mathcal{A}}}}(X)/{\Theta^{[1]}}$ to layer $0$.

\begin{definition}  We will denote by $\mathbf{T}^{(0,1)}_{\Sigma^{\boldsymbol{\mathcal{A}}}}(X)/{\Theta^{[1]}}$ the $\Sigma$-algebra 
$$
\mathbf{in}^{\Sigma,(0,1)}
\left(\mathbf{T}_{\Sigma^{\boldsymbol{\mathcal{A}}}}(X)/{\Theta^{[1]}}\right).
$$
We will call $\mathbf{T}^{(0,1)}_{\Sigma^{\boldsymbol{\mathcal{A}}}}(X)/{\Theta^{[1]}}$ the \emph{$\Sigma$-reduct} of the $\Sigma^{\boldsymbol{\mathcal{A}}}$-algebra $\mathbf{T}_{\Sigma^{\boldsymbol{\mathcal{A}}}}(X)/{\Theta^{[1]}}$. In this regard see Definition~\ref{DRed}.
\end{definition}

In the following lemmas we investigate the $\Sigma^{\boldsymbol{\mathcal{A}}}$-congruence $\Theta^{[1]}$ on 
$\mathbf{T}_{\Sigma^{\boldsymbol{\mathcal{A}}}}(X)$. In particular, we will show how this congruence can be used to describe those terms that, when mapped into $\mathrm{F}_{\Sigma^{\boldsymbol{\mathcal{A}}}}(\mathbf{Pth}_{\boldsymbol{\mathcal{A}}})$ by $\mathrm{ip}^{(1,X)@}$, return a path.

We begin with a lemma that will be fundamental later. It states that if a term is such that its image under $\mathrm{ip}^{(1,X)@}$ returns a path, then all its subterms have the same property.

\begin{restatable}{lemma}{LThetaCongSub}
\label{LThetaCongSub} 
Let $s$ be a sort in $S$ and $P$ a term in $\mathrm{T}_{\Sigma^{\boldsymbol{\mathcal{A}}}}(X)_{s}$. If $\mathrm{ip}^{(1,X)@}_{s}(P)$ is a path in $\mathrm{Pth}_{\boldsymbol{\mathcal{A}},s}$, then $\mathrm{ip}^{(1,X)@}[\mathrm{Subt}_{\Sigma^{\boldsymbol{\mathcal{A}}}}(P)]$ is a subset of 
$\mathrm{Pth}_{\boldsymbol{\mathcal{A}}}$.
\end{restatable}

\begin{proof}
We prove it by induction on the height of $P$.

\textsf{Base step of the induction.} 

If $P$ is a term of height $0$, that is, if $P\in\mathrm{B}^{0}_{\Sigma^{\boldsymbol{\mathcal{A}}}}(X)_{s}$, then
$$
\mathrm{Subt}_{\Sigma^{\boldsymbol{\mathcal{A}}}}(P)=\delta^{s,P}.
$$

That is, the unique subterm of $P$ is $P$ itself. The statement trivially holds.

\textsf{Inductive step of the induction.} 

Assume that the statement holds for terms of height up to $n$, that is, for every sort $s\in S$, if $P\in\mathrm{B}^{n}_{\Sigma^{\boldsymbol{\mathcal{A}}}}(X)_{s}$ is such that $\mathrm{ip}^{(1,X)@}_{s}(P)$ is a path then $\mathrm{ip}^{(1,X)@}[\mathrm{Subt}_{\Sigma^{\boldsymbol{\mathcal{A}}}}(P)]$ is a subset of $\mathrm{Pth}_{\boldsymbol{\mathcal{A}}}$.

Now, let $P\in\mathrm{B}^{n+1}_{\Sigma^{\boldsymbol{\mathcal{A}}}}(X)_{s}$ be a term of height $n+1$. Then there exist a unique word $\mathbf{s}\in S^{\star}-\{\lambda\}$, a unique operation symbol $\tau\in\Sigma^{\boldsymbol{\mathcal{A}}}_{\mathbf{s},s}$, and a unique family of terms $(P_{j})_{j\in\bb{\mathbf{s}}}\in\mathrm{B}^{n}_{\Sigma^{\boldsymbol{\mathcal{A}}}}(X)_{\mathbf{s}}$ such that $P=\tau^{\mathbf{T}_{\Sigma^{\boldsymbol{\mathcal{A}}}}(X)}((P_{j})_{j\in\bb{\mathbf{s}}})$. In this case, 
$$
\mathrm{Subt}_{\Sigma^{\boldsymbol{\mathcal{A}}}}(P)=
\delta^{s,P}\cup\bigcup_{j\in\bb{\mathbf{s}}}
\mathrm{Subt}_{\Sigma^{\boldsymbol{\mathcal{A}}}}(P_{j}).
$$

Note that the following chain of equivalences holds
\allowdisplaybreaks
\begin{align*}
\mathrm{ip}^{(1,X)@}_{s}(P)\in \mathrm{Pth}_{\boldsymbol{\mathcal{A}},s}
&\Leftrightarrow
\mathrm{ip}^{(1,X)@}_{s}\left(
\tau^{\mathbf{T}_{\Sigma^{\boldsymbol{\mathcal{A}}}}(X)}
\left(\left(P_{j}
\right)_{j\in\bb{\mathbf{s}}}\right)
\right)\in \mathrm{Pth}_{\boldsymbol{\mathcal{A}},s}
\tag{1}
\\&\Leftrightarrow
\tau^{\mathbf{F}_{\Sigma^{\boldsymbol{\mathcal{A}}}}\left(
\mathbf{Pth}_{\boldsymbol{\mathcal{A}}}
\right)}
\left(
\left(\mathrm{ip}^{(1,X)@}_{s_{j}}\left(
P_{j}
\right)
\right)_{j\in\bb{\mathbf{s}}}
\right)\in \mathrm{Pth}_{\boldsymbol{\mathcal{A}},s}
\tag{2}
\\&\Leftrightarrow
\left(
\mathrm{ip}^{(1,X)@}_{s_{j}}\left(
P_{j}
\right)\right)_{j\in\bb{\mathbf{s}}}
\in \mathrm{Dom}\left(
\tau^{\mathbf{Pth}_{\boldsymbol{\mathcal{A}}}}
\right).
\tag{3}
\end{align*}

The first equivalence unravels the definition of $P$; the second equation follows from the fact that, according to Definition~\ref{DIp}, $\mathrm{ip}^{(1,X)@}$ is a $\Sigma^{\boldsymbol{\mathcal{A}}}$-homomorphism; the third equivalence follows from the fact that $\tau^{\mathbf{F}_{\Sigma^{\boldsymbol{\mathcal{A}}}}(
\mathbf{Pth}_{\boldsymbol{\mathcal{A}}}
)}
(
(\mathrm{ip}^{(1,X)@}_{s_{j}}
(P_{j}
)
)_{j\in\bb{\mathbf{s}}}
)$ is a path in $\mathrm{Pth}_{\boldsymbol{\mathcal{A}},s}$ if, and only, it is the case that $(\mathrm{ip}^{(1,X)@}_{s_{j}}
(P_{j}
))_{j\in\bb{\mathbf{s}}}$ is a family of paths in the domain of $\tau^{\mathbf{Pth}_{\boldsymbol{\mathcal{A}}}}$.

This completes the proof.
\end{proof}

In the following proposition we show that $\mathrm{pr}^{\Theta^{[1]}}\circ \mathrm{CH}^{(1)}$ determines a $\Sigma^{\boldsymbol{\mathcal{A}}}$-homomorphism from $\mathbf{Pth}_{\boldsymbol{\mathcal{A}}}$ to $\mathbf{T}_{\Sigma^{\boldsymbol{\mathcal{A}}}}(X)/{\Theta^{[1]}}$.

\begin{restatable}{proposition}{PThetaCH}
\label{PThetaCH} The mapping 
\[
\mathrm{pr}^{\Theta^{[1]}}\circ \mathrm{CH}^{(1)}
\colon\mathbf{Pth}_{\boldsymbol{\mathcal{A}}}
\mor
\mathbf{T}_{\Sigma^{\boldsymbol{\mathcal{A}}}}(X)/{\Theta^{[1]}}
\]
is a $\Sigma^{\boldsymbol{\mathcal{A}}}$-homomorphism.
\end{restatable} 
\begin{proof}
We prove that this composition mapping is compatible with every operation symbol in $\Sigma^{\boldsymbol{\mathcal{A}}}$.

\textsf{The mapping $\mathrm{pr}^{\Theta^{[1]}}\circ \mathrm{CH}^{(1)}$ is a $\Sigma$-homomorphism.}

Let us note that, by Proposition~\ref{PCHHom}, $\mathrm{CH}^{(1)}$ is a $\Sigma$-homomorphism from $\mathbf{Pth}_{\boldsymbol{\mathcal{A}}}^{(0,1)}$ to $\mathbf{T}_{\Sigma^{\boldsymbol{\mathcal{A}}}}^{(0,1)}(X)$. Note also that, by Proposition~\ref{PThetaCongCatAlg}, $\mathrm{pr}^{\Theta^{[1]}}$ is a 
$\Sigma^{\boldsymbol{\mathcal{A}}}$-homomorphism. Therefore, the composition $\mathrm{pr}^{\Theta^{[1]}}\circ \mathrm{CH}^{(1)}$ is a $\Sigma$-homomorphism from $\mathbf{Pth}_{\boldsymbol{\mathcal{A}}}^{(0,1)}$ to 
$\mathbf{T}_{\Sigma^{\boldsymbol{\mathcal{A}}}}^{(0,1)}(X)/{\Theta^{[1]}}$.

\textsf{The mapping $\mathrm{pr}^{\Theta^{[1]}}\circ \mathrm{CH}^{(1)}$ is compatible with the rewrite rules.}

Let $s$ be a sort in $S$ and let $\mathfrak{p}$ be a rewrite rule in $\mathcal{A}_{s}$.
The following chain of equalities holds.
\allowdisplaybreaks
\begin{align*}
\mathrm{pr}^{\Theta^{[1]}}_{s}\left(
\mathrm{CH}^{(1)}_{s}\left(
\mathfrak{p}^{\mathbf{Pth}_{\boldsymbol{\mathcal{A}}}}
\right)
\right)&=
\mathrm{pr}^{\Theta^{[1]}}_{s}\left(
\mathfrak{p}^{\mathbf{T}_{\Sigma^{\boldsymbol{\mathcal{A}}}}(X)}
\right)
\tag{1}
\\&=
\mathfrak{p}^{\mathbf{T}_{\Sigma^{\boldsymbol{\mathcal{A}}}}(X)/{\Theta^{[1]}}}.
\tag{2}
\end{align*}

In the just stated chain of equalities, the first equality follows from Proposition~\ref{PCHA} and the second equality follows from Proposition~\ref{PThetaCongCatAlg}.

\textsf{The mapping $\mathrm{pr}^{\Theta^{[1]}}\circ \mathrm{CH}^{(1)}$ is compatible with the $0$-source.}

Let $s$ be a sort in $S$ and let us consider the $0$-source operation symbol $\mathrm{sc}^{0}_{s}$ in $\Sigma^{\boldsymbol{\mathcal{A}}}_{s,s}$. Let $\mathfrak{P}$ be a  path in $\mathrm{Pth}_{\boldsymbol{\mathcal{A}},s}$. 

The following chain of equalities holds
\allowdisplaybreaks
\begin{align*}
\mathrm{pr}^{\Theta^{[1]}}_{s}\left(
\mathrm{CH}^{(1)}_{s}\left(
\mathrm{sc}^{0\mathbf{Pth}_{\boldsymbol{\mathcal{A}}}}_{s}
\left(
\mathfrak{P}
\right)\right)
\right)&=
\mathrm{pr}^{\Theta^{[1]}}_{s}\left(
\mathrm{sc}^{0\mathbf{T}_{\Sigma^{\boldsymbol{\mathcal{A}}}}(X)}_{s}\left(
\mathrm{CH}^{(1)}_{s}\left(
\mathfrak{P}
\right)\right)
\right)
\tag{1}
\\&=
\mathrm{sc}^{0\mathbf{T}_{\Sigma^{\boldsymbol{\mathcal{A}}}}(X)/{\Theta^{[1]}}}_{s}
\left(
\mathrm{pr}^{\Theta^{[1]}}_{s}\left(
\mathrm{CH}^{(1)}_{s}\left(
\mathfrak{P}
\right)\right)
\right).
\tag{2}
\end{align*}

In the just stated chain of equalities, the first equality follows from the fact that, by Definition~\ref{DTheta}, the pair
\[
\left(\mathrm{CH}^{(1)}_{s}\left(
\mathrm{sc}^{0\mathbf{Pth}_{\boldsymbol{\mathcal{A}}}}_{s}\left(
\mathfrak{P}
\right)\right),
\mathrm{sc}^{0\mathbf{T}_{\Sigma^{\boldsymbol{\mathcal{A}}}}(X)}_{s}\left(
\mathrm{CH}^{(1)}_{s}\left(
\mathfrak{P}
\right)\right)
\right)
\]
is in $\Theta^{(1)}_{s}$; finally, the second equality follows from Proposition~\ref{PThetaCongCatAlg}.

Hence, $\mathrm{pr}^{\Theta^{[1]}}\circ \mathrm{CH}^{(1)}$ is compatible with the $0$-source operation.

\textsf{The mapping $\mathrm{pr}^{\Theta^{[1]}}\circ \mathrm{CH}^{(1)}$ is compatible with the $0$-target.}

Let $s$ be a sort in $S$ and let us consider the $0$-target operation symbol $\mathrm{tg}^{0}_{s}$ in $\Sigma^{\boldsymbol{\mathcal{A}}}_{s,s}$. Let $\mathfrak{P}$ be a path in $\mathrm{Pth}_{\boldsymbol{\mathcal{A}},s}$, then the following equality holds
\[
\mathrm{pr}^{\Theta^{[1]}}_{s}\left(
\mathrm{CH}^{(1)}_{s}\left(
\mathrm{tg}^{0\mathbf{Pth}_{\boldsymbol{\mathcal{A}}}}_{s}\left(
\mathfrak{P}
\right)\right)\right)
=
\mathrm{tg}^{0\mathbf{T}_{\Sigma^{\boldsymbol{\mathcal{A}}}}(X)/{\Theta^{[1]}}}_{s}\left(
\mathrm{pr}^{\Theta^{[1]}}_{s}\left(
\mathrm{CH}^{(1)}_{s}\left(
\mathfrak{P}
\right)\right)\right).
\]

The proof of this case is identical to that of the $0$-source.

Hence, $\mathrm{pr}^{\Theta^{[1]}}\circ \mathrm{CH}^{(1)}$ is compatible with the $0$-target operation.

\textsf{The mapping $\mathrm{pr}^{\Theta^{[1]}}\circ \mathrm{CH}^{(1)}$ is compatible with the $0$-composition.}

Let $s$ be a sort in $S$ and let us consider the $0$-composition operation symbol $\circ^{0}_{s}$ in 
$\Sigma^{\boldsymbol{\mathcal{A}}}_{ss,s}$. Let $\mathfrak{P}$ and $\mathfrak{Q}$ be two paths in 
$\mathrm{Pth}_{\boldsymbol{\mathcal{A}},s}$ such that 
\[
\mathrm{sc}^{(0,1)}_{s}\left(\mathfrak{Q}\right)
=
\mathrm{tg}^{(0,1)}_{s}\left(\mathfrak{P}\right).
\]
Then the following chain of equalities holds
\begin{flushleft}
$\mathrm{pr}^{\Theta^{[1]}}_{s}\left(
\mathrm{CH}^{(1)}_{s}\left(
\mathfrak{Q}
\circ^{0\mathbf{Pth}_{\boldsymbol{\mathcal{A}}}}_{s}
\mathfrak{P}
\right)
\right)$
\allowdisplaybreaks
\begin{align*}
\qquad&=
\mathrm{pr}^{\Theta^{[1]}}_{s}\left(
\mathrm{CH}^{(1)}_{s}\left(
\mathfrak{Q}
\right)
\circ^{0\mathbf{T}_{\Sigma^{\boldsymbol{\mathcal{A}}}}(X)}_{s}
\mathrm{CH}^{(1)}_{s}\left(
\mathfrak{P}
\right)
\right)
\tag{1}
\\&=
\mathrm{pr}^{\Theta^{[1]}}_{s}\left(
\mathrm{CH}^{(1)}_{s}\left(
\mathfrak{Q}
\right)\right)
\circ^{0\mathbf{T}_{\Sigma^{\boldsymbol{\mathcal{A}}}}(X)/{\Theta^{[1]}}}_{s}
\mathrm{pr}^{\Theta^{[1]}}_{s}\left(
\mathrm{CH}^{(1)}_{s}\left(
\mathfrak{P}
\right)\right)
.
\tag{2}
\end{align*}
\end{flushleft}

In the just stated chain of equalities, the first equality follows from the fact that, by Definition~\ref{DTheta}, the pair
\[
\left(\mathrm{CH}^{(1)}_{s}\left(
\mathfrak{Q}
\circ^{0\mathbf{Pth}_{\boldsymbol{\mathcal{A}}}}_{s}
\mathfrak{P}
\right),
\mathrm{CH}^{(1)}_{s}\left(
\mathfrak{Q}
\right)
\circ^{0\mathbf{T}_{\Sigma^{\boldsymbol{\mathcal{A}}}}(X)}_{s}
\mathrm{CH}^{(1)}_{s}\left(
\mathfrak{P}
\right)
\right)
\]
is in $\Theta^{(1)}_{s}$; finally, the second equality follows from Proposition~\ref{PThetaCongCatAlg}.

Hence, $\mathrm{pr}^{\Theta^{[1]}}\circ \mathrm{CH}^{(1)}$ is compatible with the $0$-composition operation.

All in all, we conclude that the composition
$\mathrm{pr}^{\Theta^{[1]}}\circ \mathrm{CH}^{(1)}$
is a $\Sigma^{\boldsymbol{\mathcal{A}}}$-homomorphism from $\mathbf{Pth}_{\boldsymbol{\mathcal{A}}}$ to $\mathbf{T}_{\Sigma^{\boldsymbol{\mathcal{A}}}}(X)/{\Theta^{[1]}}$.

This completes the proof.
\end{proof}

The following lemma states that if a term of sort $s$ is such that its image under $\mathrm{ip}^{(1,X)@}_{s}$ is a path, then it is related, with respect to the $\Sigma^{\boldsymbol{\mathcal{A}}}$-congruence  $\Theta^{[1]}$ on $\mathbf{T}_{\Sigma^{\boldsymbol{\mathcal{A}}}}(X)$, with a term of the same sort in $\mathrm{CH}^{(1)}[\mathrm{Pth}_{\boldsymbol{\mathcal{A}}}]$. Actually, we prove that such a term is related with its image under the action of $\mathrm{CH}^{(1)}_{s}\circ\mathrm{ip}^{(1,X)@}_{s}$.

\begin{restatable}{lemma}{LWCong}
\label{LWCong} 
Let $s$ be a sort in $S$ and $P$ a term in $\mathrm{T}_{\Sigma^{\boldsymbol{\mathcal{A}}}}(X)_{s}$. If
$
\mathrm{ip}^{(1,X)@}_{s}(P)$ is a path in $\mathrm{Pth}_{\boldsymbol{\mathcal{A}},s}
$ then
$(P,\mathrm{CH}^{(1)}_{s}(\mathrm{ip}^{(1,X)@}_{s}(P)))
\in\Theta^{[1]}_{s}.$
\end{restatable}
\begin{proof}
We prove it by induction on the height of $P$.

\textsf{Base step of the induction.}

If $z\in\mathrm{B}^{0}_{\Sigma^{\boldsymbol{\mathcal{A}}}}(X)_{s}$, then either (1) $z\in \eta^{(1,0)\sharp}[\mathrm{B}^{0}_{\Sigma}(X)]_{s}$, or (2) $z=\mathfrak{p}^{\mathbf{T}_{\Sigma^{\boldsymbol{\mathcal{A}}}}(X)}$ for some rewrite rule $\mathfrak{p}\in\mathcal{A}_{s}$. Recall that in either case, 
$\mathrm{ip}^{(1,X)@}_{s}(z)$ is a path in $\mathrm{Pth}_{\boldsymbol{\mathcal{A}},s}$.

If (1), then $z=\eta^{(1,0)\sharp}_{s}(x)$ for some term $x\in\mathrm{B}^{0}_{\Sigma}(X)_{s}$.

Hence the following chain of equalities holds
\allowdisplaybreaks
\begin{align*}
\mathrm{CH}^{(1)}_{s}\left(
\mathrm{ip}^{(1,X)@}_{s}\left(
z
\right)\right)&=
\mathrm{CH}^{(1)}_{s}\left(
\mathrm{ip}^{(1,0)\sharp}_{s}\left(x
\right)\right)
\tag{1}
\\&=\eta^{(1,0)\sharp}_{s}\left(x\right)
\tag{2}
\\&=z.
\tag{3}
\end{align*}

The first equation follows from the fact that the image of $\eta^{(1,0)\sharp}_{s}(x)$ under $\mathrm{ip}^{(1,X)@}$ is, by Proposition~\ref{PIpUZ}, given by the $(1,0)$-identity path on $x$; the second equation follows from Proposition~\ref{PCHId}; finally, the last equality recovers the definition of $z$.

By reflexivity, we conclude that 
$(P,\mathrm{CH}^{(1)}_{s}(\mathrm{ip}^{(1,X)@}_{s}(P)))
\in\Theta^{[1]}_{s}.$

If (2), then $z=\mathfrak{p}^{\mathbf{T}_{\Sigma^{\boldsymbol{\mathcal{A}}}}(X)}$, for some rewrite rule $\mathfrak{p}\in\mathcal{A}_{s}$.

Hence the following chain of equalities holds
\allowdisplaybreaks
\begin{align*}
\mathrm{CH}^{(1)}_{s}\left(
\mathrm{ip}^{(1,X)@}_{s}\left(
\mathfrak{p}^{\mathbf{T}_{\Sigma^{\boldsymbol{\mathcal{A}}}}(X)}
\right)\right)
&=\mathrm{CH}^{(1)}_{s}\left(
\mathfrak{p}^{
\mathbf{F}_{\Sigma^{\boldsymbol{\mathcal{A}}}}
\left(\mathbf{Pth}_{\boldsymbol{\mathcal{A}}}\right)
}\right)
\tag{1}
\\&
=\mathrm{CH}^{(1)}_{s}\left(
\mathfrak{p}^{\mathbf{Pth}_{\boldsymbol{\mathcal{A}}}}
\right)
\tag{2}
\\&
=\mathrm{CH}^{(1)}_{s}\left(
\mathrm{ech}^{(1,\mathcal{A})}_{s}\left(
\mathfrak{p}
\right)\right)
\tag{3}
\\&
=\mathfrak{p}^{\mathbf{T}_{\Sigma^{\boldsymbol{\mathcal{A}}}}(X)}.
\tag{4}
\end{align*}

The first equality follows from the fact that $\mathrm{ip}^{(1,X)@}$ is a $\Sigma^{\boldsymbol{\mathcal{A}}}$-homomorphism in virtue of Definition~\ref{DIp}; the second equality follows 
from the fact that the operation symbol $\mathfrak{p}$ in $\Sigma^{\boldsymbol{\mathcal{A}}}_{\lambda,s}$ is defined in the many-sorted partial $\Sigma^{\boldsymbol{\mathcal{A}}}$-algebra $\mathbf{Pth}_{\boldsymbol{\mathcal{A}}}$ according to Proposition~\ref{PPthCatAlg},
thus the interpretation of the constant operation symbol $\mathfrak{p}$ in $
\mathbf{F}_{\Sigma^{\boldsymbol{\mathcal{A}}}}
(\mathbf{Pth}_{\boldsymbol{\mathcal{A}}})$ becomes that interpretation occurring in 
$\mathbf{Pth}_{\boldsymbol{\mathcal{A}}}$;
the third equality follows from the fact that $\mathfrak{p}^{\mathbf{Pth}_{\boldsymbol{\mathcal{A}}}}$ is the echelon on $\mathfrak{p}$; finally, the last equality follows from Proposition~\ref{PCHA}.

By reflexivity, we have that 
$(P,\mathrm{CH}^{(1)}_{s}(\mathrm{ip}^{(1,X)@}_{s}(P)))
\in\Theta^{[1]}_{s}.$

\textsf{Inductive step of the induction.}

Assume that the statement holds for terms up to height $n\in \mathbb{N}$, i.e., for every  term $P$ in $\mathrm{T}_{\Sigma^{\boldsymbol{\mathcal{A}}}}(X)_{s}$ with $P\in\mathrm{B}^{n}_{\Sigma^{\boldsymbol{\mathcal{A}}}}(X)_{s}$, if
$
\mathrm{ip}^{(1,X)@}_{s}(P)$ is a path in $\mathrm{Pth}_{\boldsymbol{\mathcal{A}},s}
$ then
$(P,\mathrm{CH}^{(1)}_{s}(\mathrm{ip}^{(1,X)@}_{s}(P)))
\in\Theta^{[1]}_{s}.$

Now, we prove it for the case in which 
$P$ is a term in $\mathrm{T}_{\Sigma^{\boldsymbol{\mathcal{A}}}}(X)_{s}$ with $P\in\mathrm{B}^{n+1}_{\Sigma^{\boldsymbol{\mathcal{A}}}}(X)_{s}$ satisfying that $
\mathrm{ip}^{(1,X)@}_{s}(P)$ is a path in $\mathrm{Pth}_{\boldsymbol{\mathcal{A}},s}
$. Since $P\in\mathrm{B}^{n+1}_{\Sigma^{\boldsymbol{\mathcal{A}}}}(X)_{s}$,  then there exists a unique $\mathbf{s}\in S^{\star}-\{\lambda\}$, a unique operation symbol $\tau\in\Sigma^{\boldsymbol{\mathcal{A}}}_{\mathbf{s},s}$ and a unique family of terms $(P_{j})_{j\in\bb{\mathbf{s}}}\in\mathrm{B}^{n}_{\Sigma^{\boldsymbol{\mathcal{A}}}}(X)_{\mathbf{s}}$ satisfying that 
$$P=\tau^{\mathbf{T}_{\Sigma^{\boldsymbol{\mathcal{A}}}}(X)}
\left(\left(
P_{j}
\right)_{j\in\bb{\mathbf{s}}}
\right).$$ 

Since we are assuming that 
$
\mathrm{ip}^{(1,X)@}_{s}(P)$ is a path in $\mathrm{Pth}_{\boldsymbol{\mathcal{A}},s}
$ 
and that $P$ is equal to $\tau^{\mathbf{T}_{\Sigma^{\boldsymbol{\mathcal{A}}}}(X)}((P_{j})_{j\in\bb{\mathbf{s}}})$, we conclude, by Lemma~\ref{LThetaCongSub}, that, for every $j\in\bb{\mathbf{s}}$, $\mathrm{ip}^{(1,X)@}_{s_{j}}(P_{j})$ is a path in $\mathrm{Pth}_{\boldsymbol{\mathcal{A}},s_{j}}$ and
$$
\left(
P_{j},
\mathrm{CH}^{(1)}_{s_{j}}\left(
\mathrm{ip}^{(1,X)@}_{s_{j}}\left(
P_{j}
\right)\right)\right)
\in\Theta^{[1]}_{s_{j}}.
$$

Taking into account that $\Theta^{[1]}$ is a $\Sigma^{\boldsymbol{\mathcal{A}}}$-congruence on $\mathbf{T}_{\Sigma^{\boldsymbol{\mathcal{A}}}}(X)$, then
$$
\left(P,
\tau^{\mathbf{T}_{\Sigma^{\boldsymbol{\mathcal{A}}}}(X)}
\left(\left(
\mathrm{CH}^{(1)}_{s_{j}}\left(
\mathrm{ip}^{(1,X)@}_{s_{j}}\left(
P_{j}
\right)\right)\right)_{j\in\bb{\mathbf{s}}}\right)
\right)\in\Theta^{[1]}_{s}.
$$

On the other hand, note that the following chain of equalities hold
\begin{flushleft}
$\mathrm{pr}^{\Theta^{[1]}}_{s}\left(
\mathrm{CH}^{(1)}_{s}\left(
\mathrm{ip}^{(1,X)@}_{s}\left(
P
\right)
\right)
\right)
$
\allowdisplaybreaks
\begin{align*}
\qquad&=
\mathrm{pr}^{\Theta^{[1]}}_{s}\left(
\mathrm{CH}^{(1)}_{s}\left(
\mathrm{ip}^{(1,X)@}_{s}\left(
\tau^{\mathbf{T}_{\Sigma^{\boldsymbol{\mathcal{A}}}}(X)}
\left(
\left(
P_{j}
\right)_{j\in\bb{\mathbf{s}}}
\right)
\right)
\right)
\right)
\tag{1}
\\&=
\mathrm{pr}^{\Theta^{[1]}}_{s}\left(
\mathrm{CH}^{(1)}_{s}\left(
\tau^{\mathbf{Pth}_{\boldsymbol{\mathcal{A}}}}
\left(
\left(
\mathrm{ip}^{(1,X)@}_{s_{j}}\left(
P_{j}
\right)
\right)_{j\in\bb{\mathbf{s}}}
\right)
\right)
\right)
\tag{2}
\\&=
\tau^{\mathbf{T}_{\Sigma^{\boldsymbol{\mathcal{A}}}}(X)/{\Theta^{[1]}}}
\left(
\left(
\mathrm{pr}^{\Theta^{[1]}}_{s_{j}}\left(
\mathrm{CH}^{(1)}_{s_{j}}\left(
\mathrm{ip}^{(1,X)@}_{s_{j}}\left(
P_{j}
\right)
\right)
\right)
\right)_{j\in\bb{\mathbf{s}}}
\right)
\tag{3}
\\&=
\mathrm{pr}^{\Theta^{[1]}}_{s}\left(
\tau^{\mathbf{T}_{\Sigma^{\boldsymbol{\mathcal{A}}}}(X)}
\left(
\left(
\mathrm{CH}^{(1)}_{s_{j}}\left(
\mathrm{ip}^{(1,X)@}_{s_{j}}\left(
P_{j}
\right)
\right)
\right)_{j\in\bb{\mathbf{s}}}
\right)
\right).
\tag{4}
\end{align*}
\end{flushleft}

In the just stated chain of equalities, the first equality unravels the description of $P$; the second equality follows from the fact that, by assumption $\mathrm{ip}^{(1,X)@}_{s}(P)$ is a path and, by Definition~\ref{DIp}, $\mathrm{ip}^{(1,X)@}$ is a $\Sigma^{\boldsymbol{\mathcal{A}}}$-homomorphism; the third equality follows from the fact that, according to Proposition~\ref{PThetaCH}, $\mathrm{pr}^{\Theta^{[1]}}\circ \mathrm{CH}^{(1)}$ is a $\Sigma^{\boldsymbol{\mathcal{A}}}$-homomorphism; finally, the last equality follows from the fact that, according to Proposition~\ref{PThetaCongCatAlg}, $\mathrm{pr}^{\Theta^{[1]}}$ is a $\Sigma^{\boldsymbol{\mathcal{A}}}$-homomorphism.

Consequently, the following pair is in $\Theta^{[1]}$
\[
\left(
\mathrm{CH}^{(1)}_{s}\left(
\mathrm{ip}^{(1,X)@}_{s}\left(
P
\right)
\right)
,
\tau^{\mathbf{T}_{\Sigma^{\boldsymbol{\mathcal{A}}}}(X)}
\left(
\left(
\mathrm{CH}^{(1)}_{s_{j}}\left(
\mathrm{ip}^{(1,X)@}_{s_{j}}\left(
P_{j}
\right)
\right)
\right)_{j\in\bb{\mathbf{s}}}
\right)
\right)
\in\Theta^{[1]}_{s}.
\]

Since $\Theta^{[1]}_{s}$ is an equivalence relation, we conclude that 
\[
\left(
P,
\mathrm{CH}^{(1)}_{s}\left(
\mathrm{ip}^{(1,X)@}_{s}\left(
P
\right)
\right)
\right)\in\Theta^{[1]}_{s}.
\]

This completes the proof.
\end{proof}

In the following lemma we state that, for every sort $s\in S$, if two terms are $\Theta^{[1]}_{s}$-related and one of them, when mapped under $\mathrm{ip}^{(1,X)@}_{s}$, returns a path, then the other term has a similar behaviour. Moreover, if the just described situation happens, then these two paths will have the same image under $\mathrm{CH}^{(1)}$.

\begin{restatable}{lemma}{LThetaCong}
\label{LThetaCong} 
Let $s$ be a sort in $S$ and $P,Q$ terms in $\mathrm{T}_{\Sigma^{\boldsymbol{\mathcal{A}}}}(X)_{s}$ such that $(P,Q)\in\Theta^{[1]}_{s}$, then
\begin{itemize}
\item[(i)] $
\mathrm{ip}^{(1,X)@}_{s}(P)\in\mathrm{Pth}_{\boldsymbol{\mathcal{A}},s}
$
if, and only if, 
$
\mathrm{ip}^{(1,X)@}_{s}(Q)\in\mathrm{Pth}_{\boldsymbol{\mathcal{A}},s}
$;
\item[(ii)] If  $\mathrm{ip}^{(1,X)@}_{s}(P)$ or $\mathrm{ip}^{(1,X)@}_{s}(Q)$ is a path in $\mathrm{Pth}_{\boldsymbol{\mathcal{A}},s}$ then
$$
\mathrm{CH}^{(1)}_{s}(\mathrm{ip}^{(1,X)@}_{s}(P))=
\mathrm{CH}^{(1)}_{s}(\mathrm{ip}^{(1,X)@}_{s}(Q))
.$$
\end{itemize}
\end{restatable}

\begin{proof}
We recall from Remark~\ref{RThetaCong} that 
$$
\Theta^{[1]}=
\bigcup_{n\in\mathbb{N}}\mathrm{C}^{n}_{\Sigma^{\boldsymbol{\mathcal{A}}}}\left(
\Theta^{(1)}\right).
$$ 

We prove the statement by induction on $n\in\mathbb{N}$.

\textsf{Base step of the induction.}

Let us recall that, by Remark~\ref{RThetaCong}, we have that
$$
\mathrm{C}^{0}_{\Sigma^{\boldsymbol{\mathcal{A}}}}(\Theta^{(1)})=\Theta^{(1)}\cup\left(\Theta^{(1)}\right)^{-1}\cup\Delta_{\mathrm{T}_{\Sigma^{\boldsymbol{\mathcal{A}}}}(X)}.
$$

The statement trivially holds for a pair $(P,Q)$ in $\Delta_{\mathrm{T}_{\Sigma^{\boldsymbol{\mathcal{A}}}}(X)_{s}}$.

If the pair $(P,Q)$ is in $\Theta^{(1)}_{s}$, then, according to Definition~\ref{DTheta}, the term $P$ belongs to $\mathrm{CH}^{(1)}_{s}[\mathrm{Pth}_{\boldsymbol{\mathcal{A}},s}]$. From Proposition~\ref{PIpCH}, we know that 
$\mathrm{ip}^{(1,X)@}_{s}(P)$ is a path in $\mathrm{Pth}_{\boldsymbol{\mathcal{A}},s}$. The different cases for $Q$ are handled by Lemmas~\ref{LThetaSc},~\ref{LThetaTg} or~\ref{LThetaComp}. In any case, it follows that $\mathrm{ip}^{(1,X)@}_{s}(Q)$ is also a path in $\mathrm{Pth}_{\boldsymbol{\mathcal{A}},s}$. Moreover,
$$
\mathrm{CH}^{(1)}_{s}\left(
\mathrm{ip}^{(1,X)@}_{s}\left(
P
\right)\right)
=
\mathrm{CH}^{(1)}_{s}\left(
\mathrm{ip}^{(1,X)@}_{s}\left(
Q
\right)\right).
$$

In case $(P,Q)$ is a pair in $(\Theta^{(1)})^{-1}_{s}$ we reason in a similar way.

This completes the base case.

\textsf{Inductive step of the induction.}

Assume the statement holds for $n\in\mathbb{N}$, i.e., for every sort $s\in S$ and every pair of terms $(P,Q)$ in $\mathrm{T}_{\Sigma^{\boldsymbol{\mathcal{A}}}}(X)_{s}$ such that $(P,Q)\in\mathrm{C}^{n}_{\Sigma^{\boldsymbol{\mathcal{A}}}}(\Theta^{(1)})_{s}$ then
\begin{itemize}
\item[(i)] $
\mathrm{ip}^{(1,X)@}_{s}(P)\in\mathrm{Pth}_{\boldsymbol{\mathcal{A}},s}
$
if, and only if, 
$
\mathrm{ip}^{(1,X)@}_{s}(Q)\in\mathrm{Pth}_{\boldsymbol{\mathcal{A}},s}
$;
\item[(ii)] If  $\mathrm{ip}^{(1,X)@}_{s}(P)$ or $\mathrm{ip}^{(1,X)@}_{s}(Q)$ is a path in $\mathrm{Pth}_{\boldsymbol{\mathcal{A}},s}$ then
$$
\mathrm{CH}^{(1)}_{s}\left(\mathrm{ip}^{(1,X)@}_{s}(P)\right)=
\mathrm{CH}^{(1)}_{s}\left(\mathrm{ip}^{(1,X)@}_{s}(Q)\right)
.$$
\end{itemize}

We now prove the statement for $n+1$. Let $s$ be a sort in $S$ and let $(P,Q)$ be a pair of terms in $\mathrm{T}_{\Sigma^{\boldsymbol{\mathcal{A}}}}(X)_{s}$ such that $(P,Q)\in\mathrm{C}^{n+1}_{\Sigma^{\boldsymbol{\mathcal{A}}}}(\Theta^{(1)})_{s}$. Let us recall that, by Remark~\ref{RThetaCong}, we have that 
\begin{multline*}
\mathrm{C}^{n+1}_{\Sigma^{\boldsymbol{\mathcal{A}}}}(\Theta^{(1)})
=\left(
\mathrm{C}^{n}_{\Sigma^{\boldsymbol{\mathcal{A}}}}(\Theta^{(1)})_{s}
\circ
\mathrm{C}^{n}_{\Sigma^{\boldsymbol{\mathcal{A}}}}(\Theta^{(1)})_{s}
\right)
\cup
\\
\left(
\bigcup_{\tau\in\Sigma^{\boldsymbol{\mathcal{A}}}_{\neq\lambda,s}}
\tau^{\mathbf{T}_{\Sigma^{\boldsymbol{\mathcal{A}}}}(X)}
\times
\tau^{\mathbf{T}_{\Sigma^{\boldsymbol{\mathcal{A}}}}(X)}
\left[
\mathrm{C}^{n}_{\Sigma^{\boldsymbol{\mathcal{A}}}}(\Theta^{(1)})
_{\mathrm{ar}(\tau)}
\right]
\right)_{s\in S}.
\end{multline*}
Then either (1), $(P,Q)$ is a pair in $(
\mathrm{C}^{n}_{\Sigma^{\boldsymbol{\mathcal{A}}}}(\Theta^{(1)})_{s}
\circ
\mathrm{C}^{n}_{\Sigma^{\boldsymbol{\mathcal{A}}}}(\Theta^{(1)})_{s}
)$ or (2), $(P,Q)$ is a pair in $\tau^{\mathbf{T}_{\Sigma^{\boldsymbol{\mathcal{A}}}}(X)}
\times
\tau^{\mathbf{T}_{\Sigma^{\boldsymbol{\mathcal{A}}}}(X)}
\left[
\mathrm{C}^{n}_{\Sigma^{\boldsymbol{\mathcal{A}}}}(\Theta^{(1)})
_{\mathrm{ar}(\tau)}
\right]$ for some operation symbol $\tau\in\Sigma^{\boldsymbol{\mathcal{A}}}_{\neq\lambda,s}$.

If (1), then there exists a term $R\in \mathrm{T}_{\Sigma^{\boldsymbol{\mathcal{A}}}}(X)_{s}$ for which $(P,R)$ and $(R,Q)$ belong to $\mathrm{C}^{n}_{\Sigma^{\boldsymbol{\mathcal{A}}}}(\Theta^{(1)})_{s}$. Then by induction we have that 
$$
\mathrm{ip}^{(1,X)@}_{s}(P)\in\mathrm{Pth}_{\boldsymbol{\mathcal{A}},s}
\,\Leftrightarrow\,
\mathrm{ip}^{(1,X)@}_{s}(R)\in\mathrm{Pth}_{\boldsymbol{\mathcal{A}},s}
\,\Leftrightarrow\,
\mathrm{ip}^{(1,X)@}_{s}(Q)\in\mathrm{Pth}_{\boldsymbol{\mathcal{A}},s}
$$

Moreover, in the case one of the previous elements is a path in $\mathrm{Pth}_{\boldsymbol{\mathcal{A}},s}$, we have,  by induction, that 
$$
\mathrm{CH}^{(1)}_{s}\left(
\mathrm{ip}^{(1,X)@}_{s}\left(
P
\right)\right)=
\mathrm{CH}^{(1)}_{s}\left(
\mathrm{ip}^{(1,X)@}_{s}\left(
R
\right)\right)=
\mathrm{CH}^{(1)}_{s}\left(
\mathrm{ip}^{(1,X)@}_{s}\left(
Q
\right)\right).
$$ 

If~(2), then there exists a unique word $\mathbf{s}\in S-\{\lambda\}$, a unique operation symbol $\tau\in\Sigma_{\mathbf{s},s}$ and a unique family of pairs $((P_{j},Q_{j}))_{j\in\bb{\mathbf{s}}}$ in $\mathrm{C}^{n}_{\Sigma^{\boldsymbol{\mathcal{A}}}}(\Theta^{(1)})_{\mathbf{s}}$ for which
$$
(P,Q)=\left(
\tau^{\mathbf{T}_{\Sigma^{\boldsymbol{\mathcal{A}}}}(X)}
\left(\left(P_{j}
\right)_{j\in\bb{\mathbf{s}}}\right)
,
\tau^{\mathbf{T}_{\Sigma^{\boldsymbol{\mathcal{A}}}}(X)}
\left(\left(Q_{j}
\right)_{j\in\bb{\mathbf{s}}}\right)
\right).
$$

We will distinguish the following cases according to the different possibilities for the operation symbol $\tau\in\Sigma^{\boldsymbol{\mathcal{A}}}_{s}$. Note that either (2.1) $\tau$ is an operation symbol $\sigma\in\Sigma_{\mathbf{s},s}$; or (2.1) $\tau$ is the operation symbol of $0$-source $\mathrm{sc}^{0}_{s}\in\Sigma^{\boldsymbol{\mathcal{A}}}_{s,s}$, or (2.3) $\tau$ is the operation symbol of $0$-target $\mathrm{tg}^{0}_{s}\in\Sigma^{\boldsymbol{\mathcal{A}}}_{s,s}$, or (2.4) $\tau$ is the operation of $0$-composition $\circ^{0}_{s}\in\Sigma^{\boldsymbol{\mathcal{A}}}_{ss,s}$.

\textsf{$\tau$ is an operation symbol $\sigma\in \Sigma_{\mathbf{s},s}$.}

Let $\mathbf{s}$ be a word in $S^{\star}-\{\lambda\}$, let $\sigma$ be an operation symbol in $\Sigma_{\mathbf{s},s}$ and let $((P_{j},Q_{j}))_{j\in\bb{\mathbf{s}}}$ be a family of pairs in $\mathrm{C}_{\Sigma^{\boldsymbol{\mathcal{A}}}}^{n}(\Theta^{(1)})_{\mathbf{s}}$ for which 
$$
(P,Q)=
\left(
\sigma^{\mathbf{T}_{\Sigma^{\boldsymbol{\mathcal{A}}}}(X)}
\left(\left(
P_{j}
\right)_{j\in\bb{\mathbf{s}}}\right),
\sigma^{\mathbf{T}_{\Sigma^{\boldsymbol{\mathcal{A}}}}(X)}
\left(\left(
Q_{j}
\right)_{j\in\bb{\mathbf{s}}}\right)
\right).
$$

Note that the following chain of equivalences holds
\allowdisplaybreaks
\begin{align*}
\mathrm{ip}^{(1,X)@}_{s}(P)\in\mathrm{Pth}_{\boldsymbol{\mathcal{A}},s}
&\Leftrightarrow
\mbox{for every }j\in\bb{\mathbf{s}},\,\mathrm{ip}^{(1,X)@}_{s_{j}}(P_{j})\in\mathrm{Pth}_{\boldsymbol{\mathcal{A}},s_{j}}
\tag{1}
\\&\Leftrightarrow
\mbox{for every }j\in\bb{\mathbf{s}},\,\mathrm{ip}^{(1,X)@}_{s_{j}}(Q_{j})\in\mathrm{Pth}_{\boldsymbol{\mathcal{A}},s_{j}}
\tag{2}
\\&\Leftrightarrow
\mathrm{ip}^{(1,X)@}_{s}(Q)\in\mathrm{Pth}_{\boldsymbol{\mathcal{A}},s}.
\tag{3}
\end{align*}

In the just stated chain of equivalences, the first equivalence follows from left to right by Lemma~\ref{LThetaCongSub} and from right to left because $\sigma$ is a total operation in $\mathbf{Pth}_{\boldsymbol{\mathcal{A}}}$, according to Proposition~\ref{PPthAlg}, and $\mathrm{ip}^{(1,X)@}$ is a $\Sigma^{\boldsymbol{\mathcal{A}}}$-homomorphism, by Definition~\ref{DIp}; the second equivalence follows by induction, finally the third equivalence follows from left to right because $\sigma$ is a total operation in $\mathbf{Pth}_{\boldsymbol{\mathcal{A}}}$, according to Proposition~\ref{PPthAlg},  and $\mathrm{ip}^{(1,X)@}$ is a $\Sigma^{\boldsymbol{\mathcal{A}}}$-homomorphism, according to Definition~\ref{DIp}, and from right to left by Lemma~\ref{LThetaCongSub}.

Assume, without loss of generality, that $\mathrm{ip}^{(1,X)@}_{s}(P)$ is a  path in $\mathrm{Pth}_{\boldsymbol{\mathcal{A}},s}$. As we have seen before, this is the case exactly when, for every $j\in\bb{\mathbf{s}}$, $\mathrm{ip}^{(1,X)@}_{s_{j}}(P_{j})$ is a  path in $\mathrm{Pth}_{\boldsymbol{\mathcal{A}},s_{j}}$. By induction, we also have that, for every $j\in\bb{\mathbf{s}}$, $\mathrm{ip}^{(1,X)@}_{s_{j}}(Q_{j})$ is a  path in $\mathrm{Pth}_{\boldsymbol{\mathcal{A}},s_{j}}$. Moreover, by induction, the following equality holds
\[
\mathrm{CH}^{(1)}_{s_{j}}\left(
\mathrm{ip}^{(1,X)@}_{s_{j}}\left(
P_{j}
\right)\right)
=
\mathrm{CH}^{(1)}_{s_{j}}\left(
\mathrm{ip}^{(1,X)@}_{s_{j}}\left(
Q_{j}
\right)\right).
\]

Note that the following chain of equalities holds
\allowdisplaybreaks
\begin{align*}
\mathrm{CH}^{(1)}_{s}\left(
\mathrm{ip}^{(1,X)@}_{s}\left(
P
\right)
\right)
&=
\mathrm{CH}^{(1)}_{s}\left(
\mathrm{ip}^{(1,X)@}_{s}\left(
\sigma^{\mathbf{T}_{\Sigma^{\boldsymbol{\mathcal{A}}}}(X)}
\left(
\left(
P_{j}
\right)_{j\in\bb{\mathbf{s}}}
\right)
\right)
\right)
\tag{1}
\\&=
\mathrm{CH}^{(1)}_{s}\left(
\sigma^{\mathbf{Pth}_{\boldsymbol{\mathcal{A}}}}
\left(
\left(
\mathrm{ip}^{(1,X)@}_{s_{j}}\left(
P_{j}
\right)
\right)_{j\in\bb{\mathbf{s}}}
\right)
\right)
\tag{2}
\\&=
\mathrm{CH}^{(1)}_{s}\left(
\sigma^{\mathbf{Pth}_{\boldsymbol{\mathcal{A}}}}
\left(
\left(
\mathrm{ip}^{(1,X)@}_{s_{j}}\left(
Q_{j}
\right)
\right)_{j\in\bb{\mathbf{s}}}
\right)
\right)
\tag{3}
\\&=
\mathrm{CH}^{(1)}_{s}\left(
\mathrm{ip}^{(1,X)@}_{s}\left(
\sigma^{\mathbf{T}_{\Sigma^{\boldsymbol{\mathcal{A}}}}(X)}
\left(
\left(
Q_{j}
\right)_{j\in\bb{\mathbf{s}}}
\right)
\right)
\right)
\tag{4}
\\&=
\mathrm{CH}^{(1)}_{s}\left(
\mathrm{ip}^{(1,X)@}_{s}\left(
Q
\right)
\right).
\tag{5}
\end{align*}

In the just stated chain of equalities, the first equality unravels the description of $P$; the second equivalence follows from the fact that $\sigma$ is a total operation in $\mathbf{Pth}_{\boldsymbol{\mathcal{A}}}$, according to Proposition~\ref{PPthAlg},  and $\mathrm{ip}^{(1,X)@}$ is a $\Sigma^{\boldsymbol{\mathcal{A}}}$-homomorphism, according to Definition~\ref{DIp}, and the fact that we are assuming that $\mathrm{ip}^{(1,X)@}_{s}(P)$ is a  path in $\mathrm{Pth}_{\boldsymbol{\mathcal{A}},s}$; the third equality follows by induction and from the fact that, according to Proposition~\ref{PCHCong}, $\mathrm{Ker}(\mathrm{CH}^{(1)})$ is a $\Sigma^{\boldsymbol{\mathcal{A}}}$-congruence; the fourth equality follows from the fact that $\sigma$ is a total operation in $\mathbf{Pth}_{\boldsymbol{\mathcal{A}}}$, according to Proposition~\ref{PPthAlg},  and $\mathrm{ip}^{(1,X)@}$ is a $\Sigma^{\boldsymbol{\mathcal{A}}}$-homomorphism, according to Definition~\ref{DIp}, and the fact that we are assuming that, for every $j\in\bb{\mathbf{s}}$, $\mathrm{ip}^{(1,X)@}_{s_{j}}(Q_{j})$ is a  path in $\mathrm{Pth}_{\boldsymbol{\mathcal{A}},s_{j}}$; finally, the last equality recovers the description of $Q$.

This completes the case $\sigma\in\Sigma_{\mathbf{s},s}$.

\textsf{$\tau$ is the $0$-source operation symbol.} 

Consider the case of the $0$-source operation symbol $\mathrm{sc}^{0}_{s}$  in $\Sigma^{\boldsymbol{\mathcal{A}}}_{s,s}$. Then we can find a pair $(P',Q')$ in $\mathrm{C}_{\Sigma^{\boldsymbol{\mathcal{A}}}}^{n}(\Theta^{(1)})_{s}$ for which 
$$
(P,Q)=
\left(
\mathrm{sc}_{s}^{0\mathbf{T}_{\Sigma^{\boldsymbol{\mathcal{A}}}}(X)}\left(
P'
\right),
\mathrm{sc}_{s}^{0\mathbf{T}_{\Sigma^{\boldsymbol{\mathcal{A}}}}(X)}\left(
Q'
\right)\right).
$$

Let us note that the following chain of equivalences holds
\allowdisplaybreaks
\begin{align*}
\mathrm{ip}^{(1,X)@}_{s}(P)\in\mathrm{Pth}_{\boldsymbol{\mathcal{A}}, s}
&\Leftrightarrow
\mathrm{ip}^{(1,X)@}_{s}(P')\in\mathrm{Pth}_{\boldsymbol{\mathcal{A}}, s}
\tag{1}
\\&\Leftrightarrow
\mathrm{ip}^{(1,X)@}_{s}(Q')\in\mathrm{Pth}_{\boldsymbol{\mathcal{A}}, s}
\tag{2}
\\&\Leftrightarrow
\mathrm{ip}^{(1,X)@}_{s}(Q)\in\mathrm{Pth}_{\boldsymbol{\mathcal{A}}, s}.
\tag{3}
\end{align*}

In the just stated chain of equivalences, the first equivalence follows from left to right by Lemma~\ref{LThetaCongSub} and from right to left because $\mathrm{sc}^{0}_{s}$ is a total operation in $\mathbf{Pth}_{\boldsymbol{\mathcal{A}}}$, according to Proposition~\ref{PPthCatAlg}, and $\mathrm{ip}^{(1,X)@}$ is a $\Sigma^{\boldsymbol{\mathcal{A}}}$-homomorphism, according to Definition~\ref{DIp}; the second equivalence follows by induction, finally the third equivalence follows from left to right because $\mathrm{sc}^{0}_{s}$ is a total operation in $\mathbf{Pth}_{\boldsymbol{\mathcal{A}}}$, according to Proposition~\ref{PPthCatAlg},  and $\mathrm{ip}^{(1,X)@}$ is a $\Sigma^{\boldsymbol{\mathcal{A}}}$-homomorphism, according to Definition~\ref{DIp}, and from right to left by Lemma~\ref{LThetaCongSub}.

Assume, without loss of generality, that $\mathrm{ip}^{(1,X)@}_{s}(P)$ is a  path in $\mathrm{Pth}_{\boldsymbol{\mathcal{A}},s}$. As we have seen before, this is the case exactly when $\mathrm{ip}^{(1,X)@}_{s}(P')$ is a path in $\mathrm{Pth}_{\boldsymbol{\mathcal{A}},s}$. By induction, we also have that $\mathrm{ip}^{(1,X)@}_{s}(Q')$ is a  path in $\mathrm{Pth}_{\boldsymbol{\mathcal{A}},s}$. Moreover, by induction, the following equality holds
\[
\mathrm{CH}^{(1)}_{s}\left(
\mathrm{ip}^{(1,X)@}_{s}\left(
P'
\right)\right)
=
\mathrm{CH}^{(1)}_{s}\left(
\mathrm{ip}^{(1,X)@}_{s}\left(
Q'
\right)\right).
\]

Note that the following chain of equalities holds
\allowdisplaybreaks
\begin{align*}
\mathrm{CH}^{(1)}_{s}\left(
\mathrm{ip}^{(1,X)@}_{s}\left(
P
\right)
\right)
&=
\mathrm{CH}^{(1)}_{s}\left(
\mathrm{ip}^{(1,X)@}_{s}\left(
\mathrm{sc}^{0\mathbf{T}_{\Sigma^{\boldsymbol{\mathcal{A}}}}(X)}_{s}
\left(
P'
\right)
\right)
\right)
\tag{1}
\\&=
\mathrm{CH}^{(1)}_{s}\left(
\mathrm{sc}^{0\mathbf{Pth}_{\boldsymbol{\mathcal{A}}}}_{s}
\left(
\mathrm{ip}^{(1,X)@}_{s}\left(
P'
\right)
\right)
\right)
\tag{2}
\\&=
\mathrm{CH}^{(1)}_{s}\left(
\mathrm{sc}^{0\mathbf{Pth}_{\boldsymbol{\mathcal{A}}}}_{s}
\left(
\mathrm{ip}^{(1,X)@}_{s}\left(
Q'
\right)
\right)
\right)
\tag{3}
\\&=
\mathrm{CH}^{(1)}_{s}\left(
\mathrm{ip}^{(1,X)@}_{s}\left(
\mathrm{sc}^{0\mathbf{T}_{\Sigma^{\boldsymbol{\mathcal{A}}}}(X)}_{s}
\left(
Q'
\right)
\right)
\right)
\tag{4}
\\&=
\mathrm{CH}^{(1)}_{s}\left(
\mathrm{ip}^{(1,X)@}_{s}\left(
Q
\right)
\right).
\tag{5}
\end{align*}

In the just stated chain of equalities, the first equality unravels the description of $P$; the second equivalence follows from the fact that $\mathrm{sc}^{0}_{s}$ is a total operation in $\mathbf{Pth}_{\boldsymbol{\mathcal{A}}}$, according to Proposition~\ref{PPthCatAlg}, and $\mathrm{ip}^{(1,X)@}$ is a $\Sigma^{\boldsymbol{\mathcal{A}}}$-homomorphism, according to Definition~\ref{DIp}, and the fact that we are assuming that $\mathrm{ip}^{(1,X)@}_{s}(P)$ is a  path in $\mathrm{Pth}_{\boldsymbol{\mathcal{A}},s}$; the third equality follows by induction and from the fact that, according to Proposition~\ref{PCHCong}, $\mathrm{Ker}(\mathrm{CH}^{(1)})$ is a $\Sigma^{\boldsymbol{\mathcal{A}}}$-congruence; the fourth equality follows from the fact that $\mathrm{sc}^{0}_{s}$ is a total operation in $\mathbf{Pth}_{\boldsymbol{\mathcal{A}}}$, according to Proposition~\ref{PPthCatAlg},  and $\mathrm{ip}^{(1,X)@}$ is a $\Sigma^{\boldsymbol{\mathcal{A}}}$-homomorphism, according to Definition~\ref{DIp}, and the fact that we are assuming that $\mathrm{ip}^{(1,X)@}_{s}(Q')$ is a path in $\mathrm{Pth}_{\boldsymbol{\mathcal{A}}^{(1)},s}$; finally, the last equality recovers the description of $Q$.

This completes the case of the $0$-source.

\textsf{$\tau$ is the $0$-target operation symbol.} 

Let $\tau$ be the operation symbol $\mathrm{tg}^{0}_{s}$ in $\Sigma^{\boldsymbol{\mathcal{A}}}_{s,s}$ and let $(P',Q')$ be a family of pairs in $\mathrm{C}^{n}_{\Sigma^{\boldsymbol{\mathcal{A}}}}(\Theta^{(1)})_{s}$ for which 
$$
(P,Q)=
\left(
\mathrm{tg}^{0\mathbf{T}_{\Sigma^{\boldsymbol{\mathcal{A}}}}(X)}_{s}\left(
P'
\right)
,
\mathrm{tg}^{0\mathbf{T}_{\Sigma^{\boldsymbol{\mathcal{A}}}}(X)}_{s}\left(
Q'
\right)
\right).
$$

Then, the following properties hold
\begin{itemize}
\item[(i)] $\mathrm{ip}^{(1,X)@}_{s}(P)\in\mathrm{Pth}_{\boldsymbol{\mathcal{A}},s}$ if, and only if, $\mathrm{ip}^{(1,X)@}_{s}(Q)\in\mathrm{Pth}_{\boldsymbol{\mathcal{A}},s}$;
\item[(ii)] If $\mathrm{ip}^{(1,X)@}_{s}(P)\in\mathrm{Pth}_{\boldsymbol{\mathcal{A}},s}$ or  $\mathrm{ip}^{(1,X)@}_{s}(Q)\in\mathrm{Pth}_{\boldsymbol{\mathcal{A}},s}$ is a  path in $\mathrm{Pth}_{\boldsymbol{\mathcal{A}}^{(1)},s}$ then
$
\mathrm{CH}^{(1)}_{s}(
\mathrm{ip}^{(1,X)@}_{s}(
P))=\mathrm{CH}^{(1)}_{s}(
\mathrm{ip}^{(1,X)@}_{s}(
Q)).
$
\end{itemize}

The proof of this case is similar to that of the $0$-source.

This completes the case of the $0$-target.

\textsf{$\tau$ is the $0$-composition operation symbol.} 

Finally, consider the case of the $0$-composition operation symbol $\circ^{0}_{s}$ in $\Sigma^{\boldsymbol{\mathcal{A}}}_{ss,s}$. Then we can find two pairs $(P',Q')$ and $(P',Q')$ in $\mathrm{C}_{\Sigma^{\boldsymbol{\mathcal{A}}}}^{n}(\Theta^{(1)})_{s}$ for which 
$$
(P,Q)=
\left(P''\circ_{s}^{0\mathbf{T}_{\Sigma^{\boldsymbol{\mathcal{A}}}}(X)}P',
Q''\circ_{s}^{0\mathbf{T}_{\Sigma^{\boldsymbol{\mathcal{A}}}}(X)}Q'
\right).
$$

Let us note that the following chain of equivalences holds
\allowdisplaybreaks
\begin{align*}
\mathrm{ip}^{(1,X)@}_{s}(P)
\in
\mathrm{Pth}_{\mathcal{A},s}
&\Leftrightarrow
\begin{cases}
\mathrm{ip}^{(1,X)@}_{s}(P''),
\mathrm{ip}^{(1,X)@}_{s}(P')
\in
\mathrm{Pth}_{\mathcal{A},s}
\\
\mathrm{sc}^{(0,1)}_{s}(\mathrm{ip}^{(1,X)@}_{s}(P''))=
\mathrm{tg}^{(0,1)}_{s}(\mathrm{ip}^{(1,X)@}_{s}(P'))
\end{cases}
\tag{1}
\\&\Leftrightarrow
\begin{cases}
\mathrm{ip}^{(1,X)@}_{s}(Q''),
\mathrm{ip}^{(1,X)@}_{s}(Q')
\in
\mathrm{Pth}_{\mathcal{A},s}
\\
\mathrm{sc}^{(0,1)}_{s}(\mathrm{ip}^{(1,X)@}_{s}(Q''))=
\mathrm{tg}^{(0,1)}_{s}(\mathrm{ip}^{(1,X)@}_{s}(Q'))
\end{cases}
\tag{2}
\\&\Leftrightarrow
\mathrm{ip}^{(1,X)@}_{s}(Q)
\in
\mathrm{Pth}_{\mathcal{A},s}.
\tag{3}
\end{align*}

In the just stated chain of equivalences, the first equivalence follows from left to right by Lemma~\ref{LThetaCongSub} and by the description of the $0$-composition operation in $\mathbf{Pth}_{\boldsymbol{\mathcal{A}}}$, according to Proposition~\ref{PPthCatAlg}, and from right to left because $\circ^{0}_{s}$ is a defined operation in $\mathbf{Pth}_{\boldsymbol{\mathcal{A}}}$, according to Proposition~\ref{PPthCatAlg}, and $\mathrm{ip}^{(1,X)@}$ is a $\Sigma^{\boldsymbol{\mathcal{A}}}$-homomorphism, according to Definition~\ref{DIp}; the second equivalence follows by induction. Note also that
\allowdisplaybreaks
\begin{align*}
\mathrm{CH}^{(1)}_{s}\left(
\mathrm{ip}^{(1,X)@}_{s}\left(
P'
\right)
\right)
&=
\mathrm{CH}^{(1)}_{s}\left(
\mathrm{ip}^{(1,X)@}_{s}\left(
Q'
\right)
\right);
\\
\mathrm{CH}^{(1)}_{s}\left(
\mathrm{ip}^{(1,X)@}_{s}\left(
P''
\right)
\right)
&=
\mathrm{CH}^{(1)}_{s}\left(
\mathrm{ip}^{(1,X)@}_{s}\left(
Q''
\right)
\right).
\end{align*}
Therefore, according to Lemma~\ref{LCH}, we have that 
\allowdisplaybreaks
\begin{align*}
\mathrm{sc}^{(0,1)}_{s}\left(
\mathrm{ip}^{(1,X)@}_{s}\left(
P''
\right)\right)&=
\mathrm{sc}^{(0,1)}_{s}\left(
\mathrm{ip}^{(1,X)@}_{s}\left(
Q''
\right)\right)
\\
\mathrm{tg}^{(0,1)}_{s}\left(
\mathrm{ip}^{(1,X)@}_{s}\left(
P'
\right)\right)&=
\mathrm{tg}^{(0,1)}_{s}\left(
\mathrm{ip}^{(1,X)@}_{s}\left(
Q'
\right)\right);
\end{align*}
finally the third equivalence follows from left to right because $\circ^{0}_{s}$ is a defined operation in $\mathbf{Pth}_{\boldsymbol{\mathcal{A}}}$, according to Proposition~\ref{PPthCatAlg}, and $\mathrm{ip}^{(1,X)@}$ is a $\Sigma^{\boldsymbol{\mathcal{A}}}$-homomorphism, according to Definition~\ref{DIp}, and from right to left by Lemma~\ref{LThetaCongSub} and by the description of the $0$-composition operation in $\mathbf{Pth}_{\boldsymbol{\mathcal{A}}}$, according to Proposition~\ref{PPthCatAlg}.

Assume, without loss of generality, that $\mathrm{ip}^{(1,X)@}_{s}(P)$ is a path in $\mathrm{Pth}_{\boldsymbol{\mathcal{A}},s}$. As we have seen before, this is the case exactly when $\mathrm{ip}^{(1,X)@}_{s}(P')$ and $\mathrm{ip}^{(1,X)@}_{s}(P'')$ are  paths in $\mathrm{Pth}_{\boldsymbol{\mathcal{A}},s}$ satisfying that 
\[
\mathrm{sc}^{(0,1)}_{s}\left(
\mathrm{ip}^{(1,X)@}_{s}\left(
P''
\right)\right)=
\mathrm{tg}^{(0,1)}_{s}\left(
\mathrm{ip}^{(1,X)@}_{s}\left(
P'
\right)\right).
\]
By induction, we also have that $\mathrm{ip}^{(1,X)@}_{s}(Q')$ and $\mathrm{ip}^{(1,X)@}_{s}(Q'')$ are paths in $\mathrm{Pth}_{\boldsymbol{\mathcal{A}},s}$. As we have proven before, their $0$-composition is well-defined. Moreover, by induction, the following equality holds
\allowdisplaybreaks
\begin{align*}
\mathrm{CH}^{(1)}_{s}\left(
\mathrm{ip}^{(1,X)@}_{s}\left(
P'
\right)\right)
&=
\mathrm{CH}^{(1)}_{s}\left(
\mathrm{ip}^{(1,X)@}_{s}\left(
Q'
\right)\right);
\\
\mathrm{CH}^{(1)}_{s}\left(
\mathrm{ip}^{(1,X)@}_{s}\left(
P''
\right)\right)
&=
\mathrm{CH}^{(1)}_{s}\left(
\mathrm{ip}^{(1,X)@}_{s}\left(
Q''
\right)\right);
\end{align*}

Note that the following chain of equalities holds
\allowdisplaybreaks
\begin{align*}
\mathrm{CH}^{(1)}_{s}\left(
\mathrm{ip}^{(1,X)@}_{s}\left(
P
\right)
\right)
&=
\mathrm{CH}^{(1)}_{s}\left(
\mathrm{ip}^{(1,X)@}_{s}\left(
P''
\circ^{0\mathbf{T}_{\Sigma^{\boldsymbol{\mathcal{A}}}}(X)}_{s}
P'
\right)
\right)
\tag{1}
\\&=
\mathrm{CH}^{(1)}_{s}\left(
\mathrm{ip}^{(1,X)@}_{s}\left(
P''
\right)
\circ^{0\mathbf{Pth}_{\boldsymbol{\mathcal{A}}}}_{s}
\mathrm{ip}^{(1,X)@}_{s}\left(
P'
\right)
\right)
\tag{2}
\\&=
\mathrm{CH}^{(1)}_{s}\left(
\mathrm{ip}^{(1,X)@}_{s}\left(
Q''
\right)
\circ^{0\mathbf{Pth}_{\boldsymbol{\mathcal{A}}}}_{s}
\mathrm{ip}^{(1,X)@}_{s}\left(
Q'
\right)
\right)
\tag{3}
\\&=
\mathrm{CH}^{(1)}_{s}\left(
\mathrm{ip}^{(1,X)@}_{s}\left(
Q''
\circ^{0\mathbf{T}_{\Sigma^{\boldsymbol{\mathcal{A}}}}(X)}_{s}
Q'
\right)
\right)
\tag{4}
\\&=
\mathrm{CH}^{(1)}_{s}\left(
\mathrm{ip}^{(1,X)@}_{s}\left(
Q
\right)
\right).
\tag{5}
\end{align*}

In the just stated chain of equalities, the first equality unravels the description of $P$; the second equivalence follows from the fact that $\circ^{0}_{s}$ is a defined operation in $\mathbf{Pth}_{\boldsymbol{\mathcal{A}}}$, according to Proposition~\ref{PPthCatAlg},  and $\mathrm{ip}^{(1,X)@}$ is a $\Sigma^{\boldsymbol{\mathcal{A}}}$-homomorphism, according to Definition~\ref{DIp}, and the fact that we are assuming that $\mathrm{ip}^{(1,X)@}_{s}(P)$ is a  path in $\mathrm{Pth}_{\boldsymbol{\mathcal{A}},s}$; the third equality follows by induction and from the fact that, according to Proposition~\ref{PCHCong}, $\mathrm{Ker}(\mathrm{CH}^{(1)})$ is a $\Sigma^{\boldsymbol{\mathcal{A}}}$-congruence; the fourth equality follows from the fact that $\circ^{0}_{s}$ is a defined operation in $\mathbf{Pth}_{\boldsymbol{\mathcal{A}}}$, according to Proposition~\ref{PPthCatAlg},  and $\mathrm{ip}^{(1,X)@}$ is a $\Sigma^{\boldsymbol{\mathcal{A}}}$-homomorphism, according to Definition~\ref{DIp}, and the fact that we are assuming that $\mathrm{ip}^{(1,X)@}_{s}(Q')$ and $\mathrm{ip}^{(1,X)@}_{s}(Q'')$ are paths in $\mathrm{Pth}_{\boldsymbol{\mathcal{A}},s}$ whose $0$-composition is well-defined; finally, the last equality recovers the description of $Q$.

This completes the case of the $0$-composition.

This concludes the proof.
\end{proof}

\begin{restatable}{corollary}{CThetaCong}
\label{CThetaCong} Let $s$ be a sort in $S$, $\mathfrak{P}$ a path in $\mathrm{Pth}_{\boldsymbol{\mathcal{A}},s}$ and $P$ a term in $\mathrm{T}_{\Sigma^{\boldsymbol{\mathcal{A}}}}(X)_{s}$ such that 
$
(
P,
\mathrm{CH}^{(1)}_{s}(\mathfrak{P})
)
\in 
\Theta^{[1]}_{s}.
$
Then $\mathrm{ip}^{(1,X)@}_{s}(P)$ is a path in $[\mathfrak{P}]_{s}$.
\end{restatable}
\begin{proof}
By Proposition~\ref{PIpCH}, $\mathrm{ip}^{(1,X)@}_{s}(\mathrm{CH}^{(1)}_{s}(\mathfrak{P}))$ is a path in $[\mathfrak{P}]_{s}$. By Lemma~\ref{LThetaCong}, we have that $\mathrm{ip}^{(1,X)@}_{s}(P)$ is a path in $\mathrm{Pth}_{\boldsymbol{\mathcal{A}},s}$. Moreover, the following chain of equalities holds
$$
\mathrm{CH}^{(1)}_{s}\left(
\mathrm{ip}^{(1,X)@}_{s}\left(
P
\right)\right)
=
\mathrm{CH}^{(1)}_{s}\left(
\mathrm{ip}^{(1,X)@}_{s}\left(
\mathrm{CH}^{(1)}_{s}\left(
\mathfrak{P}
\right)\right)\right)
=
\mathrm{CH}^{(1)}_{s}\left(
\mathfrak{P}
\right).
$$
This completes the proof.
\end{proof}

\begin{restatable}{corollary}{CThetaCongCH}
\label{CThetaCongCH} 
Let $s$ be a sort in $S$ and $\mathfrak{P}',\mathfrak{P}$ two paths in $\mathrm{Pth}_{\boldsymbol{\mathcal{A}},s}$. If the pair
$
(
\mathrm{CH}^{(1)}_{s}(\mathfrak{P}'),
\mathrm{CH}^{(1)}_{s}(\mathfrak{P})
)$ is in $
\Theta^{[1]}_{s}
$, then $\mathrm{CH}^{(1)}_{s}(\mathfrak{P}')=\mathrm{CH}^{(1)}_{s}(\mathfrak{P})$.
\end{restatable}

\begin{proof}
The following chain of equalities holds
\allowdisplaybreaks
\begin{align*}
\mathrm{CH}^{(1)}_{s}\left(
\mathfrak{P}'
\right)&=
\mathrm{CH}^{(1)}_{s}\left(
\mathrm{ip}^{(1,X)@}_{s}\left(
\mathrm{CH}^{(1)}_{s}\left(
\mathfrak{P}'
\right)\right)\right)
\tag{1}
\\&
=\mathrm{CH}^{(1)}_{s}\left(
\mathfrak{P}\right).
\tag{2}
\end{align*}

The first equation follows from Proposition~\ref{PIpCH}, note that $\mathrm{ip}^{(1,X)@}_{s}(\mathrm{CH}^{(1)}_{s}(\mathfrak{P}))$  is a path in $[\mathfrak{P}']^{}_{s}$; the second equation follows from Corollary~\ref{CThetaCong}.
\end{proof}

\begin{remark}\label{RExplFalla} 
From Lemmas~\ref{LWCong} and~\ref{LThetaCong} and Corollary~\ref{CThetaCongCH} we infer that, for every sort $s\in S$, if a term $P\in\mathrm{T}_{\Sigma^{\boldsymbol{\mathcal{A}}}}(X)_{s}$ belongs to $[\mathrm{CH}^{(1)}_{s}(\mathfrak{P})]_{\Theta^{[1]}_{s}}$, for some path $\mathfrak{P}\in\mathrm{Pth}_{\boldsymbol{\mathcal{A}},s}$, then $P$ can be understood as an alternative term describing the same path as $\mathfrak{P}$ up to equivalence with respect to $\mathrm{Ker}(\mathrm{CH}^{(1)})$.
Therefore, we can think of the equivalence class $[\mathrm{CH}^{(1)}_{s}(\mathfrak{P})]_{\Theta^{[1]}_{s}}$ as a set containing all possible terms describing all paths in $[\mathfrak{P}]_{s}$. Thus, for a term $P\in\mathrm{T}_{\Sigma^{\boldsymbol{\mathcal{A}}}}(X)_{s}$ such that $P\in [\mathrm{CH}^{(1)}_{s}(\mathfrak{P})]_{\Theta^{[1]}_{s}}$, we have that $\mathrm{ip}^{(1,X)@}_{s}(P)$ is a path in $[\mathfrak{P}]_{s}$. 
\end{remark}			
\chapter{
\texorpdfstring
{The $S$-sorted set $\mathrm{PT}_{\boldsymbol{\mathcal{A}}}$ of path terms}
{Path terms}
}\label{S1I}

In this chapter we introduce the $S$-sorted set of path terms, denoted by  $\mathrm{PT}_{\boldsymbol{\mathcal{A}}}$, as the $\Theta^{[1]}$-saturation of the subset $\mathrm{CH}^{(1)}[\mathrm{Pth}_{\boldsymbol{\mathcal{A}}}]$ of $\mathrm{T}_{\Sigma^{\boldsymbol{\mathcal{A}}}}(X)$---these terms will play a crucial role in the iteration of the results presented here. We show that a term $P$ in $\mathrm{T}_{\Sigma^{\boldsymbol{\mathcal{A}}}}(X)_{s}$ is a path term if, and only if, $\mathrm{ip}^{(1,X)@}_{s}(P)$ is a path in $\mathrm{Pth}_{\boldsymbol{\mathcal{A}},s}$. In fact this allows us to consider the corestriction of $\mathrm{ip}^{(1,X)@}$ to $\mathrm{Pth}_{\boldsymbol{\mathcal{A}}}$ and the correstriction of $\mathrm{CH}^{(1)}$ to $\mathrm{PT}_{\boldsymbol{\mathcal{A}}}$, that will keep the same notation. We prove that $\mathrm{PT}_{\boldsymbol{\mathcal{A}}}$ has an structure of many-sorted partial $\Sigma^{\boldsymbol{\mathcal{A}}}$-algebra, denoted by $\mathbf{PT}_{\boldsymbol{\mathcal{A}}}$, inherited from $\mathbf{T}_{\Sigma^{\boldsymbol{\mathcal{A}}}}(X)/{\Theta^{[1]}}$. We can prove that $\mathbf{PT}_{\boldsymbol{\mathcal{A}}}$ is a partial Dedekind-Peano $\Sigma^{\boldsymbol{\mathcal{A}}}$-algebra. This implies, in particular, that subterms of path terms are itself path terms. It also implies that the many-sorted set $\mathrm{PT}_{\boldsymbol{\mathcal{A}}}$ of path terms also admits an Artinian partial order, denoted by $\leq_{\mathbf{PT}_{\boldsymbol{\mathcal{A}}}}$, given by the standard subterm relation. We then introduce the quotient of path terms by the restriction of the relation $\Theta^{[1]}$, that we will denote by $[\mathrm{PT}_{\boldsymbol{\mathcal{A}}}]$. The possible compositions of $\mathrm{ip}^{(1,X)@}$ and $\mathrm{CH}^{(1)}$ with $\mathrm{pr}^{\Theta^{[1]}}$, the projection to the quotient, will be denoted by $\mathrm{ip}^{([1],X)@}$ and $\mathrm{CH}^{[1]}$, respectively. It is shown that $[\mathrm{PT}_{\boldsymbol{\mathcal{A}}}]$ has a structure of many-sorted partial $\Sigma^{\boldsymbol{\mathcal{A}}}$-algebra, denoted $[\mathbf{PT}_{\boldsymbol{\mathcal{A}}}]$. Moreover, $[\mathrm{PT}_{\boldsymbol{\mathcal{A}}}]$ is an $S$-sorted category, more precisely, it has a structure of categorial $\Sigma$-algebra, denoted by $[\mathsf{PT}_{\boldsymbol{\mathcal{A}}}]$. Finally, the quotient $[\mathrm{PT}_{\boldsymbol{\mathcal{A}}}]$ admits an Artinian partial order, denoted by $\leq_{[\mathbf{PT}_{\boldsymbol{\mathcal{A}}}]}$, defined using representatives and the partial order  $\leq_{\mathbf{PT}_{\boldsymbol{\mathcal{A}}}}$ introduced before.

We next introduce the notion of path term. For a sort $s\in S$, a term in $\mathrm{T}_{\Sigma^{\boldsymbol{\mathcal{A}}}}(X)_{s}$ will be called a path term if it is a term that can be interpreted as a path in $\mathrm{Pth}_{\boldsymbol{\mathcal{A}},s}$ by means of $\mathrm{ip}^{(1,X)@}$.

\begin{restatable}{definition}{DPT}
\label{DPT} 
\index{path terms!first-order!$\mathrm{PT}_{\boldsymbol{\mathcal{A}}}$}
We let $\mathrm{PT}_{\boldsymbol{\mathcal{A}}}$ stand for $\left[\mathrm{CH}^{(1)}\left[
\mathrm{Pth}_{\boldsymbol{\mathcal{A}}}\right]\right]^{\Theta^{[1]}}$, the $\Theta^{[1]}$-saturation of the subset $\mathrm{CH}^{(1)}[\mathrm{Pth}_{\boldsymbol{\mathcal{A}}}]$ of $\mathrm{T}_{\Sigma^{\boldsymbol{\mathcal{A}}}}(X)$,  and we call it the $S$-sorted set of \emph{path terms}.
Thus, for a sort $s\in S$, a term $P\in\mathrm{T}_{\Sigma^{\boldsymbol{\mathcal{A}}}}(X)_{s}$ belongs to  
$\mathrm{PT}_{\boldsymbol{\mathcal{A}},s}$ if, and only if, there exists a path $\mathfrak{P}\in\mathrm{Pth}_{\boldsymbol{\mathcal{A}},s}$ such that 
$$\left(
P,\mathrm{CH}^{(1)}_{s}\left(
\mathfrak{P}
\right)\right)\in\Theta^{[1]}_{s}.$$
\end{restatable}

In the following proposition we prove the most important property of path terms. For every sort $s\in S$, a term $P$ in $\mathrm{T}_{\Sigma^{\boldsymbol{\mathcal{A}}}}(X)_{s}$ is a path term if, and only if, $\mathrm{ip}^{(1,X)@}_{s}(P)$ is a path.

\begin{restatable}{proposition}{PPT}
\label{PPT} 
Let $s$ be a sort in $S$ and $P$ a term in $\mathrm{T}_{\Sigma^{\boldsymbol{\mathcal{A}}}}(X)_{s}$. Then the following conditions are equivalent.
\begin{enumerate}
\item[(i)] $P\in\mathrm{PT}_{\boldsymbol{\mathcal{A}},s}$.
\item[(ii)] $\mathrm{ip}^{(1,X)@}_{s}(P)$ is a path in $\mathrm{Pth}_{\boldsymbol{\mathcal{A}},s}$.
\end{enumerate} 
\end{restatable}

\begin{proof}
Assume that $P$ is a term in $\mathrm{PT}_{\boldsymbol{\mathcal{A}},s}$. Then there exists a path $\mathfrak{P}\in\mathrm{Pth}_{\boldsymbol{\mathcal{A}},s}$ for which $(P,\mathrm{CH}^{(1)}_{s}(\mathfrak{P}))$ is a pair in $\Theta^{[1]}_{s}$. Then, by Corollary~\ref{CThetaCong}, $\mathrm{ip}^{(1,X)@}_{s}(P)$ is a path in $\mathrm{Pth}_{\boldsymbol{\mathcal{A}},s}$. Now, assume that $\mathrm{ip}^{(1,X)@}_{s}(P)$ is a path in $\mathrm{Pth}_{\boldsymbol{\mathcal{A}},s}$. Then, by Lemma~\ref{LWCong}, $(P,\mathrm{CH}^{(1)}_{s}(\mathrm{ip}^{(1,X)@}_{s}(P)))$ is a pair in $\Theta^{[1]}_{s}$, i.e., $P$ is a term in $\mathrm{PT}_{\boldsymbol{\mathcal{A}},s}$.
\end{proof}

In the following corollary we state that two path terms in the same $\Theta^{[1]}$-class, when interpreted as paths by means of $\mathrm{ip}^{(1,X)@}$, have the same $(0,1)$-source and  $(0,1)$-target.

\begin{restatable}{corollary}{CPTScTg}
\label{CPTScTg}
Let $s$ be a sort in $S$ and $Q$, $P$ path terms in $\mathrm{PT}_{\boldsymbol{\mathcal{A}},s}$ in the same 
$\Theta^{[1]}$-equivalence class, that is, 
$(Q,P)\in\Theta^{[1]}_{s}$. Then 
$$
\left(\mathrm{ip}^{(1,X)@}_{s}\left(
Q
\right),
\mathrm{ip}^{(1,X)@}_{s}\left(
P
\right)\right)
\in\mathrm{Ker}\left(
\mathrm{sc}^{(0,1)}
\right)_{s}
\cap\mathrm{Ker}\left(
\mathrm{tg}^{(0,1)}
\right)_{s}.
$$
\end{restatable}
\begin{proof}
By Proposition~\ref{PPT} and Lemma~\ref{LThetaCong}, it is the case that 
$$
\mathrm{CH}^{(1)}_{s}\left(
\mathrm{ip}^{(1,X)@}_{s}\left(
Q
\right)\right)
=
\mathrm{CH}^{(1)}_{s}\left(
\mathrm{ip}^{(1,X)@}_{s}\left(
P
\right)\right).
$$
This implies, by Lemma~\ref{LCH}, that 
$$
\left(
\mathrm{ip}^{(1,X)@}_{s}\left(
Q
\right),\mathrm{ip}^{(1,X)@}_{s}\left(
P
\right)\right)
\in\mathrm{Ker}\left(
\mathrm{sc}^{(0,1)}
\right)_{s}
\cap\mathrm{Ker}\left(
\mathrm{tg}^{(0,1)}
\right)_{s}.
$$
This completes the proof.
\end{proof}

In the following corollaries we state that some already known $S$-sorted mappings have nice restrictions, corestrictions or birestrictions to the $S$-sorted set of path terms. We will start with the embeddings introduced in Definition~\ref{DEta}.

\begin{restatable}{corollary}{CPTEta}
\label{CPTEta} The $S$-sorted embeddings
$\eta^{(1,X)}$ and $\eta^{(1,\mathcal{A})}$ introduced in Definition~\ref{DEta} from, respectively, $X$ and $\mathcal{A}$ to $\mathrm{T}_{\Sigma^{\boldsymbol{\mathcal{A}}}}(X)$ corestrict to $\mathrm{PT}_{\boldsymbol{\mathcal{A}}}$.
\end{restatable}
\begin{proof}
The following equations hold by Propositions~\ref{PCHId} and~\ref{PCHA}, respectively
\begin{align*}
\eta^{(1,X)}&=
\mathrm{CH}^{(1)}\circ\mathrm{ip}^{(1,X)};
&
\eta^{(1,\mathcal{A})}&=
\mathrm{CH}^{(1)}\circ\mathrm{ech}^{(1,\mathcal{A})}.
\end{align*}
Therefore, $\eta^{(1,X)}$ and $\eta^{(1,\mathcal{A})}$ always return terms in $\mathrm{CH}^{(1)}[\mathrm{Pth}_{\boldsymbol{\mathcal{A}}}]$.
\end{proof}

On the basis of the above corollary, we can consider the corestrictions of the embeddings introduced in Definition~\ref{DEta} to path terms. Let us note that, although we have a formal mechanism to denote the corestriction of a mapping to a subset of its codomain, we have decided to keep the existing notations for such corestrictions, since the labels to be used are complex and the context, so we hope, will avoid any confusion (the same considerations apply to the restrictions and birestrictions).

\begin{restatable}{definition}{DPTEta}
\label{DPTEta}  Let $X$ be an $S$-sorted set and let $\mathrm{PT}_{\boldsymbol{\mathcal{A}}}$ be the $S$-sorted set of path terms.  From now on, we will denote by
\begin{enumerate}
\item[(i)] $\eta^{(1,X)}$ the $S$-sorted mapping from $X$ to 
$\mathrm{PT}_{\boldsymbol{\mathcal{A}}}$ that, for every sort $s\in S$, sends an element $x\in X_{s}$ to the variable  $x\in \mathrm{PT}_{\boldsymbol{\mathcal{A}},s}$.
\index{identity!first-order!$\eta^{(1,X)}$}
\item[(ii)] $\eta^{(1,\mathcal{A})}$ the $S$-sorted mapping from 
$\mathcal{A}$ to 
$\mathrm{PT}_{\boldsymbol{\mathcal{A}}}$ 
that, for every sort $s\in S$, sends a  rewrite rule $\mathfrak{p}\in \mathcal{A}_{s}$ to the constant $\mathfrak{p}^{\mathrm{T}_{\Sigma^{\boldsymbol{\mathcal{A}}}}(X)}$.
\index{inclusion!first-order!$\eta^{(1,\mathcal{A})}$}
\end{enumerate}
The above $S$-sorted mappings are depicted in the diagram of Figure~\ref{FPTEta}.
\end{restatable}

\begin{figure}
\begin{tikzpicture}
[ACliment/.style={-{To [angle'=45, length=5.75pt, width=4pt, round]}},scale=.8]
\node[] (x) at (0,0) [] {$X$};
\node[] (a) at (0,-1.5) [] {$\mathcal{A}$};
\node[] (T) at (6,-1.5) [] {$\mathrm{PT}_{\boldsymbol{\mathcal{A}}}$};
\draw[ACliment, bend left=10]  (x) to node [above right] {$\eta^{(1,X)}$} (T);
\draw[ACliment]  (a) to node [below] {$\eta^{(1,\mathcal{A})}$} (T);
\end{tikzpicture}
\caption{Refined embeddings relative to $X$ and $\mathcal{A}$ at layer 1.}
\label{FPTEta}
\end{figure}

We next consider the embedding introduced in Proposition~\ref{PEmb}.

\begin{restatable}{corollary}{CPTEtas}
\label{CPTEtas} The $S$-sorted embedding $\eta^{(1,0)\sharp}$ introduced in Proposition~\ref{PEmb} from $\mathrm{T}_{\Sigma}(X)$ to $\mathrm{T}_{\Sigma^{\boldsymbol{\mathcal{A}}}}(X)$ corestricts to $\mathrm{PT}_{\boldsymbol{\mathcal{A}}}$.
\end{restatable}
\begin{proof}
The following equality holds by Proposition~\ref{PIpUZ}
$$
\mathrm{ip}^{(1,X)@}\circ\eta^{(1,0)\sharp}
=
\eta^{(1,\mathbf{Pth}_{\boldsymbol{\mathcal{A}}})}
\circ
\mathrm{ip}^{(1,0)\sharp}.
$$

Therefore, for every sort $s\in S$, and every term $P\in\mathrm{T}_{\Sigma}(X)_{s}$, we have that 
$\mathrm{ip}^{(1,X)@}_{s}(\eta^{(1,0)\sharp}_{s}(P))$ is a path in $\mathrm{Pth}_{\boldsymbol{\mathcal{A}},s}$. Hence, by Proposition~\ref{PPT}, we have that $\eta^{(1,0)\sharp}_{s}(P)$ is a path term in $\mathrm{PT}_{\boldsymbol{\mathcal{A}},s}$.

This completes the proof.
\end{proof}

By the above corollary, we can consider the corestriction of the embedding introduced in Proposition~\ref{PEmb} to path terms. To avoid further notation, we resignify the already existing many-sorted mapping.

\begin{restatable}{definition}{DPTEtas}
\label{DPTEtas} Let $X$ be an $S$-sorted set and let $\mathrm{PT}_{\boldsymbol{\mathcal{A}}}$ be the many-sorted set of path terms. From now on, we will denote by
\begin{enumerate}
\item $\eta^{(1,0)\sharp}$ the $S$-sorted mapping from $\mathrm{T}_{\Sigma}(X)$ to $\mathrm{PT}_{\boldsymbol{\mathcal{A}}}$ that, for every sort $s\in S$, sends a term $P\in\mathrm{T}_{\Sigma}(X)_{s}$ to the path term $\eta^{(1,0)\sharp}_{s}(P)$ in $\mathrm{PT}_{\boldsymbol{\mathcal{A}},s}$.
\index{inclusion!first-order!$\eta^{(1,0)\sharp}$}
\end{enumerate}
\end{restatable}

We next consider the case corresponding to $\mathrm{ip}^{(1,X)@}$.

\begin{restatable}{corollary}{CPTIp}
\label{CPTIp} 
The restriction of $\mathrm{ip}^{(1,X)@}$ to the set of path terms, that is,
$$
\mathrm{ip}^{(1,X)@}{\,\!\upharpoonright}_{\mathrm{PT}_{\boldsymbol{\mathcal{A}}}}
\colon
\mathrm{PT}_{\boldsymbol{\mathcal{A}}}
\mor
\mathrm{F}_{\Sigma^{\boldsymbol{\mathcal{A}}}}(\mathbf{Pth}_{\boldsymbol{\mathcal{A}}}),
$$
corestricts to $\mathrm{Pth}_{\boldsymbol{\mathcal{A}}}$.
\end{restatable}
\begin{proof}
It follows from Proposition~\ref{PPT}.
\end{proof}

\begin{definition}\label{DPTIp}
\index{identity!first-order!$\mathrm{ip}^{(1,X)a}$}
We agreed to also use $\mathrm{ip}^{(1,X)@}$ for the birestriction of the $S$-sorted mapping 
$\mathrm{ip}^{(1,X)@}$ to $\mathrm{PT}_{\boldsymbol{\mathcal{A}}}$ and $\mathrm{Pth}_{\boldsymbol{\mathcal{A}}}$. Thus 
$$
\mathrm{ip}^{(1,X)@}
\colon
\mathrm{PT}_{\boldsymbol{\mathcal{A}}}
\mor
\mathrm{Pth}_{\boldsymbol{\mathcal{A}}}.
$$
\end{definition}

Now we turn to $\mathrm{CH}^{(1)}$.

\begin{restatable}{corollary}{CPTCH}
\label{CPTCH}
The Curry-Howard mapping, that is,
$$
\mathrm{CH}^{(1)}\colon\mathrm{Pth}_{\boldsymbol{\mathcal{A}}}
\mor
\mathrm{T}_{\Sigma^{\boldsymbol{\mathcal{A}}}}(X),
$$
corestricts to $\mathrm{PT}_{\boldsymbol{\mathcal{A}}}$.
\end{restatable}
\begin{proof}
It suffices to take into account that $\mathrm{CH}^{(1)}[\mathrm{Pth}_{\boldsymbol{\mathcal{A}}}]$ is a subset of its own saturation with respect to the $\Sigma^{\boldsymbol{\mathcal{A}}}$-congruence $\Theta^{[1]}$.
\end{proof}

\begin{definition}\label{DPTCH}
\index{Curry-Howard!first-order!$\mathrm{CH}^{(1)}$}
We agreed to also use $\mathrm{CH}^{(1)}$  for the corestriction of the $S$-sorted mapping $\mathrm{CH}^{(1)}$ to $\mathrm{PT}_{\boldsymbol{\mathcal{A}}}$. Thus 
$$
\mathrm{CH}^{(1)}\colon\mathrm{Pth}_{\boldsymbol{\mathcal{A}}}
\mor
\mathrm{PT}_{\boldsymbol{\mathcal{A}}}.
$$

\end{definition}

\begin{figure}
\begin{tikzpicture}
[ACliment/.style={-{To [angle'=45, length=5.75pt, width=4pt, round]}},scale=0.8]
\node[] (x) at (0,0) [] {$X$};
\node[] (t) at (6,0) [] {$\mathrm{T}_{\Sigma}(X)$};
\node[] (pth1) at (6,-2) [] {$\mathrm{Pth}_{\boldsymbol{\mathcal{A}}}
$};
\node[] (pt1) at (6,-4) [] {$\mathrm{PT}_{\boldsymbol{\mathcal{A}}}$};

\draw[ACliment]  (x) to node [above right]
{$\textstyle \eta^{(0,X)}$} (t);
\draw[ACliment, bend right=10]  (x) to node [pos=.5, fill=white]
{$\textstyle \mathrm{ip}^{(1,X)}$} (pth1.west);
\draw[ACliment, bend right=20]  (x) to node  [below left, fill=white]
{$\textstyle \eta^{(1,X)}$} (pt1.west);

\draw[ACliment] 
($(t)+(0,-.35)$) to node [right] {$\textstyle \mathrm{ip}^{(1,0)\sharp}$}  ($(pth1)+(0,.35)$);

\draw[ACliment] 
($(pth1)+(-.30,-.35)$) to node [left] {$\textstyle \mathrm{CH}^{(1)}$}  ($(pt1)+(-.30,.35)$);
\draw[ACliment] 
($(pt1)+(.30,+.35)$) to node [right] {$\textstyle \mathrm{ip}^{(1,X)@}$}  ($(pth1)+(.30,-.35)$);

\draw[ACliment, rounded corners]
(t) -- 
($(t)+(2,0)$) -- node [right] {$\textstyle \eta^{(1,0)\sharp}$} 
($(pt1)+(2,0)$) -- ($(pt1.east)+(0,0)$)
;

\end{tikzpicture}
\caption{Refined embeddings relative to $X$ at layers $0$ \& $1$.}
\label{FPTX}
\end{figure}

\begin{figure}
\begin{tikzpicture}
[ACliment/.style={-{To [angle'=45, length=5.75pt, width=4pt, round]}},scale=0.8]
\node[] (a) at (0,-2) [] {$\mathcal{A}$};
\node[] (pth1) at (6,-2) [] {$\mathrm{Pth}_{\boldsymbol{\mathcal{A}}}
$};
\node[] (pt1) at (6,-4) [] {$\mathrm{PT}_{\boldsymbol{\mathcal{A}}}$};

\draw[ACliment]  (a) to node [above]
{$\textstyle \mathrm{ech}^{(1,\mathcal{A})}$} (pth1);
\draw[ACliment, bend right=10]  (a) to node [below left]
{$\textstyle \eta^{(1,\mathcal{A})}$} (pt1.west);

\draw[ACliment] 
($(pth1)+(-.30,-.35)$) to node [left] {$\textstyle \mathrm{CH}^{(1)}$}  ($(pt1)+(-.30,.35)$);
\draw[ACliment] 
($(pt1)+(.30,+.35)$) to node [right] {$\textstyle \mathrm{ip}^{(1,X)@}$}  ($(pth1)+(.30,-.35)$);

\end{tikzpicture}
\caption{Refined embeddings relative to $\mathcal{A}$ at layer $1$.}
\label{FPTA}
\end{figure}

\begin{remark} The above refinements are considered in Figures~\ref{FPTX} and~\ref{FPTA}, relatively to $X$ and $\mathcal{A}$, respectively. Every possible triangle considered in these figures commute.
\end{remark}

\section{
\texorpdfstring
{A structure of partial $\Sigma^{\boldsymbol{\mathcal{A}}}$-algebra on $\mathrm{PT}_{\boldsymbol{\mathcal{A}}}$}
{A partial algebra on path terms}
}

We next show that $\mathrm{PT}_{\boldsymbol{\mathcal{A}}}$ is equipped with a structure of partial $\Sigma^{\boldsymbol{\mathcal{A}}}$-algebra, which is a subalgebra of $\mathbf{T}_{\Sigma^{\boldsymbol{\mathcal{A}}}}(X)$. 

\begin{restatable}{proposition}{PPTCatAlg}
\index{path terms!first-order!$\mathbf{PT}_{\boldsymbol{\mathcal{A}}}$}
\label{PPTCatAlg} The $S$-sorted set $\mathrm{PT}_{\boldsymbol{\mathcal{A}}}$ is equipped, in a natural way, with a structure of partial $\Sigma^{\boldsymbol{\mathcal{A}}}$-algebra, which is a $\Sigma^{\boldsymbol{\mathcal{A}}}$-subalgebra of $\mathbf{T}_{\Sigma^{\boldsymbol{\mathcal{A}}}}(X)$.
\end{restatable}
\begin{proof}
Let us denote by $\mathbf{PT}_{\boldsymbol{\mathcal{A}}}$ the partial $\Sigma^{\boldsymbol{\mathcal{A}}}$-algebra defined as follows:

\textsf{(i)} The underlying $S$-sorted set of $\mathbf{PT}_{\boldsymbol{\mathcal{A}}}$ is $\mathrm{PT}_{\boldsymbol{\mathcal{A}}}$.

\textsf{(ii)} For every $(\mathbf{s},s)\in S^{\star}\times S$ and every operation symbol $\sigma\in\Sigma_{\mathbf{s},s}$, the operation $\sigma^{\mathbf{PT}_{\boldsymbol{\mathcal{A}}}}$ is equal to $\sigma^{\mathbf{T}_{\Sigma^{\boldsymbol{\mathcal{A}}}}(X)}$, i.e., to the interpretation of $\sigma$ in $\mathbf{T}_{\Sigma^{\boldsymbol{\mathcal{A}}}}(X)$. 

Let us note that, if $(P_{j})_{j\in\bb{\mathbf{s}}}$ is a family of path terms in $\mathrm{PT}_{\boldsymbol{\mathcal{A}},\mathbf{s}}$, then the following chain of equalities holds
\allowdisplaybreaks
\begin{align*}
\mathrm{ip}^{(1,X)@}_{s}\left(
\sigma^{\mathbf{T}_{\Sigma^{\boldsymbol{\mathcal{A}}}}(X)}
\left(\left(
P_{j}
\right)_{j\in\bb{\mathbf{s}}}
\right)\right)
&=
\sigma^{\mathbf{F}_{\Sigma^{\boldsymbol{\mathcal{A}}}}
\left(\mathbf{Pth}_{\boldsymbol{\mathcal{A}}}\right)}
\left(\left(
\mathrm{ip}^{(1,X)@}_{s_{j}}\left(
P_{j}
\right)\right)_{j\in\bb{\mathbf{s}}}\right)
\tag{1}
\\&=
\sigma^{\mathbf{Pth}_{\boldsymbol{\mathcal{A}}}}
\left(\left(
\mathrm{ip}^{(1,X)@}_{s_{j}}\left(
P_{j}
\right)\right)_{j\in\bb{\mathbf{s}}}\right).
\tag{2}
\end{align*}

The first equality follows from the fact that, by Definition~\ref{DIp}, $\mathrm{ip}^{(1,X)@}$ is a $\Sigma^{\boldsymbol{\mathcal{A}}}$-homomorphism; finally, the second equality  follows from the fact that, for every $j\in\bb{\mathbf{s}}$, the term $P_{j}$ is a path term in $\mathrm{PT}_{\boldsymbol{\mathcal{A}},s_{j}}$. Hence, in virtue of Proposition~\ref{PPT}, for every $j\in\bb{\mathbf{s}}$, $\mathrm{ip}^{(1,X)@}_{s_{j}}(P_{j})$ is a path in $\mathrm{Pth}_{\boldsymbol{\mathcal{A}},s_{j}}$. Thus, the interpretation of the operation symbol $\sigma$ in $\mathbf{F}_{\Sigma^{\boldsymbol{\mathcal{A}}}}(\mathbf{Pth}_{\boldsymbol{\mathcal{A}}})$ is given by the interpretation of $\sigma$ in $\mathbf{Pth}_{\boldsymbol{\mathcal{A}}}$.

Hence, since $\mathrm{ip}^{(1,X)@}_{s}(\sigma^{\mathbf{T}_{\Sigma^{\boldsymbol{\mathcal{A}}}}(X)}((P_{j})_{j\in\bb{\mathbf{s}}}))$ is a path in $\mathrm{Pth}_{\boldsymbol{\mathcal{A}},s}$, we conclude, following Proposition~\ref{PPT}, that $\sigma^{\mathbf{T}_{\Sigma^{\boldsymbol{\mathcal{A}}}}(X)}((P_{j})_{j\in\bb{\mathbf{s}}})$ is a path term in $\mathrm{PT}_{\boldsymbol{\mathcal{A}},s}$. Thus, the operation $\sigma^{\mathbf{PT}_{\boldsymbol{\mathcal{A}}}}$ is well-defined.

\textsf{(iii)} For every $s\in S$ and every rewrite rule $\mathfrak{p}\in \mathcal{A}_{s}$, the constant $\mathfrak{p}^{\mathbf{PT}_{\boldsymbol{\mathcal{A}}}}$ is equal to $\mathfrak{p}^{\mathbf{T}_{\Sigma^{\boldsymbol{\mathcal{A}}}}(X)}$, i.e., to the interpretation of $\mathfrak{p}$ in $\mathbf{T}_{\Sigma^{\boldsymbol{\mathcal{A}}}}(X)$. 

Let us note that the following chain of equalities holds
\allowdisplaybreaks
\begin{align*}
\mathrm{ip}^{(1,X)@}_{s}\left(
\mathfrak{p}^{\mathbf{T}_{\Sigma^{\boldsymbol{\mathcal{A}}}}(X)}
\right)
&=
\mathfrak{p}^{\mathbf{F}_{\Sigma^{\boldsymbol{\mathcal{A}}}}
\left(\mathbf{Pth}_{\boldsymbol{\mathcal{A}}}\right)}
\tag{1}
\\&=
\mathfrak{p}^{\mathbf{Pth}_{\boldsymbol{\mathcal{A}}}}.
\tag{2}
\end{align*}

The first equality follows from the fact that, by Definition~\ref{DIp}, $\mathrm{ip}^{(1,X)@}$ is a $\Sigma^{\boldsymbol{\mathcal{A}}}$-homomorphism; finally, the second equality follows from the fact that the interpretation of the constant operation symbol $\mathfrak{p}$ in $\mathbf{F}_{\Sigma^{\boldsymbol{\mathcal{A}}}}(\mathbf{Pth}_{\boldsymbol{\mathcal{A}}})$ is given by the interpretation of $\mathfrak{p}$ in $\mathbf{Pth}_{\boldsymbol{\mathcal{A}}}$.

Thus, since $\mathrm{ip}^{(1,X)@}_{s}(\mathfrak{p}^{\mathbf{T}_{\Sigma^{\boldsymbol{\mathcal{A}}}}(X)})$ is a path in $\mathrm{Pth}_{\boldsymbol{\mathcal{A}},s}$, we conclude, by Proposition~\ref{PPT}, that $\mathfrak{p}^{\mathbf{T}_{\Sigma^{\boldsymbol{\mathcal{A}}}}(X)}$ is a path term in $\mathrm{PT}_{\boldsymbol{\mathcal{A}},s}$. That is, the constant $\mathfrak{p}^{\mathbf{T}_{\Sigma^{\boldsymbol{\mathcal{A}}}}(X)}$ is well-defined.

\textsf{(iv)} For every $s\in S$, the operation $\mathrm{sc}_{s}^{0\mathbf{PT}_{\boldsymbol{\mathcal{A}}}}$ is equal to $\mathrm{sc}_{s}^{0\mathbf{T}_{\Sigma^{\boldsymbol{\mathcal{A}}}}(X)}$, i.e., to the interpretation of $\mathrm{sc}_{s}^{0}$ in $\mathbf{T}_{\Sigma^{\boldsymbol{\mathcal{A}}}}(X)$. 

Let us note that, if $P$ is a path term in $\mathrm{PT}_{\boldsymbol{\mathcal{A}},s}$, then the following chain of equalities holds
\allowdisplaybreaks
\begin{align*}
\mathrm{ip}^{(1,X)@}_{s}\left(
\mathrm{sc}_{s}^{0\mathbf{T}_{\Sigma^{\boldsymbol{\mathcal{A}}}}(X)}\left(
P
\right)\right)
&=
\mathrm{sc}_{s}^{0\mathbf{F}_{\Sigma^{\boldsymbol{\mathcal{A}}}}
\left(\mathbf{Pth}_{\boldsymbol{\mathcal{A}}}\right)}\left(
\mathrm{ip}^{(1,X)@}_{s}\left(
P
\right)\right)
\tag{1}
\\&=
\mathrm{sc}_{s}^{0\mathbf{Pth}_{\boldsymbol{\mathcal{A}}}}\left(
\mathrm{ip}^{(1,X)@}_{s}\left(
P
\right)\right).
\tag{2}
\end{align*}

The first equality follows from the fact that, by Definition~\ref{DIp}, $\mathrm{ip}^{(1,X)@}$ is a $\Sigma^{\boldsymbol{\mathcal{A}}}$-homomorphism; finally, the second equality follows from the fact that $P$ is a path term in $\mathrm{PT}_{\boldsymbol{\mathcal{A}},s}$. Hence, in virtue of Proposition~\ref{PPT}, $\mathrm{ip}^{(1,X)@}_{s}(P)$ is a path in $\mathrm{Pth}_{\boldsymbol{\mathcal{A}},s}$. thus, the interpretation of the operation symbol $\mathrm{sc}^{0}_{s}$ in $\mathbf{F}_{\Sigma^{\boldsymbol{\mathcal{A}}}}(\mathbf{Pth}_{\boldsymbol{\mathcal{A}}})$ is given by the interpretation of $\mathrm{sc}^{0}_{s}$ in $\mathbf{Pth}_{\boldsymbol{\mathcal{A}}}$.

Therefore, since $\mathrm{ip}^{(1,X)@}_{s}(\mathrm{sc}_{s}^{0\mathbf{T}_{\Sigma^{\boldsymbol{\mathcal{A}}}}(X)}(P))$ is a path in $\mathrm{Pth}_{\boldsymbol{\mathcal{A}},s}$, we have that, by Proposition~\ref{PPT},  $\mathrm{sc}_{s}^{0\mathbf{T}_{\Sigma^{\boldsymbol{\mathcal{A}}}}(X)}(P)$ is a path term in $\mathrm{PT}_{\boldsymbol{\mathcal{A}},s}$. That is, the operation $\mathrm{sc}_{s}^{0\mathbf{PT}_{\boldsymbol{\mathcal{A}}}}$ is well-defined.

\textsf{(v)} For every $s\in S$, the operation $\mathrm{tg}_{s}^{0\mathbf{PT}_{\boldsymbol{\mathcal{A}}}}$ is equal to $\mathrm{tg}_{s}^{0\mathbf{T}_{\Sigma^{\boldsymbol{\mathcal{A}}}}(X)}$, i.e., to the interpretation of $\mathrm{tg}_{s}^{0}$ in $\mathbf{T}_{\Sigma^{\boldsymbol{\mathcal{A}}}}(X)$.  The operation $\mathrm{tg}_{s}^{0\mathbf{PT}_{\boldsymbol{\mathcal{A}}}}$ is well-defined by a similar argument to that of $\mathrm{sc}_{s}^{0\mathbf{PT}_{\boldsymbol{\mathcal{A}}}}$.

\textsf{(vi)} For every $s\in S$, the partial binary operation $\circ_{s}^{0\mathbf{PT}_{\boldsymbol{\mathcal{A}}}}$ is defined on path terms $Q,P\in\mathrm{PT}_{\boldsymbol{\mathcal{A}},s}$ if, and only if, the following equation holds
$$
\mathrm{sc}_{s}^{(0,1)}\left(
\mathrm{ip}^{(1,X)@}_{s}\left(
Q
\right)\right)=
\mathrm{tg}_{s}^{(0,1)}\left(
\mathrm{ip}^{(1,X)@}_{s}\left(
P
\right)\right).
$$
For this case, $\circ_{s}^{0\mathbf{PT}_{\boldsymbol{\mathcal{A}}}}$ is equal to $\circ_{s}^{\mathbf{T}_{\Sigma^{\boldsymbol{\mathcal{A}}}}(X)}$, i.e., to the interpretation of $\circ^{0}_{s}$ in $\mathbf{T}_{\Sigma^{\boldsymbol{\mathcal{A}}}}(X)$. 

Let us note that, if $Q$ and $P$ are path terms in $\mathrm{PT}_{\boldsymbol{\mathcal{A}},s}$ with 
$$\mathrm{sc}_{s}^{(0,1)}\left(
\mathrm{ip}^{(1,X)@}_{s}\left(
Q
\right)\right)=
\mathrm{tg}_{s}^{(0,1)}\left(
\mathrm{ip}^{(1,X)@}_{s}\left(
P
\right)\right),$$
then the following chain of equalities holds
\allowdisplaybreaks
\begin{align*}
\mathrm{ip}^{(1,X)@}_{s}\left(
Q\circ_{s}^{0\mathbf{T}_{\Sigma^{\boldsymbol{\mathcal{A}}}}(X)} P
\right)
&=
\left(\mathrm{ip}^{(1,X)@}_{s}\left(
Q
\right)\right)
\circ_{s}^{0\mathbf{F}_{\Sigma^{\boldsymbol{\mathcal{A}}}}
\left(\mathbf{Pth}_{\boldsymbol{\mathcal{A}}}\right)}
\left(\mathrm{ip}^{(1,X)@}_{s}\left(
P
\right)\right)
\tag{1}
\\&=
\left(\mathrm{ip}^{(1,X)@}_{s}\left(
Q
\right)\right)
\circ_{s}^{0\mathbf{Pth}_{\boldsymbol{\mathcal{A}}}}
\left(\mathrm{ip}^{(1,X)@}_{s}\left(
P
\right)\right).
\tag{2}
\end{align*}

The first equality follows from the fact that, by Definition~\ref{DIp}, $\mathrm{ip}^{(1,X)@}$ is a $\Sigma^{\boldsymbol{\mathcal{A}}}$-homomorphism; finally, the second equality follows from the fact that $Q$ and $P$ are path terms in $\mathrm{PT}_{\boldsymbol{\mathcal{A}},s}$. Hence, in virtue of Proposition~\ref{PPT}, $\mathrm{ip}^{(1,X)@}_{s}(Q)$ and $\mathrm{ip}^{(1,X)@}_{s}(P)$ are paths in $\mathrm{Pth}_{\boldsymbol{\mathcal{A}},s}$. Moreover, $\mathrm{sc}_{s}^{(0,1)}(\mathrm{ip}^{(1,X)@}_{s}(Q))=
\mathrm{tg}_{s}^{(0,1)}(\mathrm{ip}^{(1,X)@}_{s}(P))$, i.e., these paths can be $0$-composed. Thus, the interpretation of the operation symbol $\circ^{0}_{s}$ in $\mathbf{F}_{\Sigma^{\boldsymbol{\mathcal{A}}}}(\mathbf{Pth}_{\boldsymbol{\mathcal{A}}})$ is given by the interpretation of $\circ^{0}_{s}$ in $\mathbf{Pth}_{\boldsymbol{\mathcal{A}}}$.

Since $\mathrm{ip}^{(1,X)@}_{s}(Q\circ_{s}^{0\mathbf{T}_{\Sigma^{\boldsymbol{\mathcal{A}}}}(X)} P)$ is a path in $\mathrm{Pth}_{\boldsymbol{\mathcal{A}},s}$, we have, by Proposition~\ref{PPT}, that $Q\circ_{s}^{0\mathbf{T}_{\Sigma^{\boldsymbol{\mathcal{A}}}}(X)} P$ is a path term in $\mathrm{PT}_{\boldsymbol{\mathcal{A}},s}$. That is, the operation $\circ_{s}^{0\mathbf{PT}_{\boldsymbol{\mathcal{A}}}}$ is well-defined.

This completes the definition of the partial $\Sigma^{\boldsymbol{\mathcal{A}}}$-algebra $\mathbf{PT}_{\boldsymbol{\mathcal{A}}}$. Thus defined it is obvious that $\mathbf{PT}_{\boldsymbol{\mathcal{A}}}$ is a 
$\Sigma^{\boldsymbol{\mathcal{A}}}$-subalgebra of $\mathbf{T}_{\Sigma^{\boldsymbol{\mathcal{A}}}}(X)$.
\end{proof}

We next prove that the partial $\Sigma^{\boldsymbol{\mathcal{A}}}$-algebra $\mathbf{PT}_{\boldsymbol{\mathcal{A}}}$ is a partial Dedekind-Peano $\Sigma^{\boldsymbol{\mathcal{A}}}$-algebra.

\begin{restatable}{proposition}{PPTPDP}
\label{PPTPDP} 
The partial $\Sigma^{\boldsymbol{\mathcal{A}}}$-algebra $\mathbf{PT}_{\boldsymbol{\mathcal{A}}}$ is a partial Dedekind-Peano $\Sigma^{\boldsymbol{\mathcal{A}}}$-algebra.
\end{restatable}

\begin{proof} 
According to Definition~\ref{DPAlgDP}, we need to check that the following items hold

\textsf{PDP1.} For every $(\mathbf{s},s)\in S^{\star}\times S$ and every $\tau\in\Sigma^{\boldsymbol{\mathcal{A}}}_{\mathbf{s},s}$, the operation 
$$\tau^{\mathbf{PT}_{\boldsymbol{\mathcal{A}}}}\colon
\mathbf{PT}_{\boldsymbol{\mathcal{A}},\mathbf{s}}
\dmor
\mathbf{PT}_{\boldsymbol{\mathcal{A}},s}
$$
is injective, i.e., for every $(Q_{j})_{j\in\bb{\mathbf{s}}}, (P_{j})_{j\in\bb{\mathbf{s}}}$ in $\mathrm{Dom}(\tau)$, if 
$$\tau^{\mathbf{PT}_{\boldsymbol{\mathcal{A}}}}
\left(\left(
Q_{j}
\right)_{j\in\bb{\mathbf{s}}}\right)
=\tau^{\mathbf{PT}_{\boldsymbol{\mathcal{A}}}}
\left(\left(
P_{j}
\right)_{j\in\bb{\mathbf{s}}}\right),$$
then $(Q_{j})_{j\in\bb{\mathbf{s}}}=(P_{j})_{j\in\bb{\mathbf{s}}}$. 

\textsf{PDP2.} For every $s\in S$ and every $\sigma,\tau\in\Sigma^{\boldsymbol{\mathcal{A}}}_{\cdot, s}$, if $\sigma\neq \tau$, then 
$$
\sigma^{\mathbf{PT}_{\boldsymbol{\mathcal{A}}}}\left[
\mathrm{Dom}\left(
\sigma^{\mathbf{PT}_{\boldsymbol{\mathcal{A}}}}
\right)\right]
\cap
\tau^{\mathbf{PT}_{\boldsymbol{\mathcal{A}}}}\left[
\mathrm{Dom}\left(
\tau^{\mathbf{PT}_{\boldsymbol{\mathcal{A}}}}
\right)\right]
=
\varnothing.
$$

\textsf{PDP3.} The following equality holds
$$
\textstyle
\mathrm{Sg}_{\mathbf{PT}_{\boldsymbol{\mathcal{A}}}}
\left(
\mathrm{PT}_{\boldsymbol{\mathcal{A}}}
-
\left(\bigcup_{\tau\in\Sigma_{\cdot,s}}
\mathrm{Im}\left(
\tau^{\mathbf{PT}_{\boldsymbol{\mathcal{A}}}}
\right)
\right)_{s\in S}
\right)
=
\mathrm{PT}_{\boldsymbol{\mathcal{A}}}.
$$
This is true, because this happens for the free $\Sigma^{\boldsymbol{\mathcal{A}}}$-algebra 
$\mathbf{T}_{\Sigma^{\boldsymbol{\mathcal{A}}}}(X)$. Note that the basis of the Dedekind-Peano of $\mathbf{PT}_{\boldsymbol{\mathcal{A}}}$ is 
$$
\textstyle
\mathrm{PT}_{\boldsymbol{\mathcal{A}}}
-
\left(\bigcup_{\tau\in\Sigma_{\cdot,s}}
\mathrm{Im}\left(
\tau^{\mathbf{PT}_{\boldsymbol{\mathcal{A}}}}
\right)
\right)_{s\in S}=X.
$$
\end{proof}

\section{
\texorpdfstring
{An Artinian order on $\coprod\mathrm{PT}_{\boldsymbol{\mathcal{A}}}$}
{An Artinian order on path terms}
}

In this section we define an Artinian order on $\coprod\mathrm{PT}_{\boldsymbol{\mathcal{A}}}$. To begin with we state in the following corollary that the subterms of a  path term are also path terms.

\begin{restatable}{corollary}{CPTSubt}
\label{CPTSubt} 
Let $s$ be a sort in $S$ and $P$ a path term in $\mathrm{PT}_{\boldsymbol{\mathcal{A}},s}$. Then 
$\mathrm{Subt}_{\Sigma^{\boldsymbol{\mathcal{A}}}}(P)\subseteq \mathrm{PT}_{\boldsymbol{\mathcal{A}}}$.
\end{restatable}
\begin{proof}
It follows from Definition~\ref{DPT} and Lemma~\ref{LThetaCongSub}.
\end{proof}

Consequently, we can consider the restriction of the Artinian order $\leq_{\mathbf{T}_{\Sigma^{\boldsymbol{\mathcal{A}}}}(X)}$ on $\coprod\mathrm{T}_{\Sigma^{\boldsymbol{\mathcal{A}}}}(X)$ to $\coprod\mathrm{PT}_{\boldsymbol{\mathcal{A}}}$.

\begin{restatable}{definition}{DPTOrd}
\label{DPTOrd}
\index{partial order!first-order!$\leq_{\mathbf{PT}_{\boldsymbol{\mathcal{A}}}}$}
Let $\leq_{\mathbf{PT}_{\boldsymbol{\mathcal{A}}}}$ be the binary relation on $\mathrm{PT}_{\boldsymbol{\mathcal{A}}}$ containing every pair $((Q,t),(P,s))$ in $(\coprod\mathrm{PT}_{\boldsymbol{\mathcal{A}}})^{2}$ such that
$$
\left(
Q,t
\right)
\leq_{\mathbf{T}_{\Sigma^{\boldsymbol{\mathcal{A}}}}(X)}
\left(
P,s
\right).
$$
Thus, $(Q,t)$ $\leq_{\mathbf{PT}_{\boldsymbol{\mathcal{A}}}}$-precedes $(P,s)$ if, and only if, $Q$ is a subterm of type $t$ of $P$.
\end{restatable}

\begin{restatable}{proposition}{PPTOrdArt}
\label{PPTOrdArt}
$(\coprod\mathrm{PT}_{\boldsymbol{\mathcal{A}}}, \leq_{\mathbf{PT}_{\boldsymbol{\mathcal{A}}}})$ is an ordered set. Moreover, in it there is not any strictly decreasing $\omega_{0}$-chain, i.e., $(\coprod\mathrm{PT}_{\boldsymbol{\mathcal{A}}}, \leq_{\mathbf{PT}_{\boldsymbol{\mathcal{A}}}})$ is an Artinian ordered set.
\end{restatable}

\section{
\texorpdfstring
{On the quotient $[\mathrm{PT}_{\boldsymbol{\mathcal{A}}}]$}
{On the quotient of path terms}
}

In this section we define the $S$-sorted set of path term classes as the quotient of 
$\mathrm{PT}_{\boldsymbol{\mathcal{A}}}$ by the restriction of $\Theta^{[1]}$, the 
$\Sigma^{\boldsymbol{\mathcal{A}}}$-congruence on $\mathbf{T}_{\Sigma^{\boldsymbol{\mathcal{A}}}}(X)$ introduced in Definition~\ref{DThetaCong}, to it. In addition, we investigate the $S$-sorted mappings from and to the just mentioned quotient of 
$\mathrm{PT}_{\boldsymbol{\mathcal{A}}}$.

\begin{restatable}{convention}{CPTClass}
\label{CPTClass} 
\index{path-term!first-order!$[P]_{s}$}
\index{path-term!first-order!$[\mathrm{PT}_{\boldsymbol{\mathcal{A}}}]$}
To simplify the notation, for a sort $s\in S$ and a path term $P\in\mathrm{PT}_{\boldsymbol{\mathcal{A}},s}$, we will let $[P]_{s}$ stand for $[P]_{\Theta^{[1]}_{s}}$, the 
$\Theta^{[1]}_{s}$-equivalence class of $P$, and we will call it the \emph{path term class of} $P$. 
\end{restatable}
 
\begin{restatable}{definition}{DPTQuot}
\label{DPTQuot} 
\index{path terms!first-order!$[\mathrm{PT}_{\boldsymbol{\mathcal{A}}}]$}
We denote by $[\mathrm{PT}_{\boldsymbol{\mathcal{A}}}]$ the image of 
$\mathrm{PT}_{\boldsymbol{\mathcal{A}}}$ under $\mathrm{pr}^{\Theta^{[1]}}$, the canonical projection from 
$\mathrm{T}_{\Sigma^{\boldsymbol{\mathcal{A}}}}(X)$ to $\mathrm{T}_{\Sigma^{\boldsymbol{\mathcal{A}}}}(X)/\Theta^{[1]}$,  i.e., $[\mathrm{PT}_{\boldsymbol{\mathcal{A}}}]=\mathrm{pr}^{\Theta^{[1]}}[\mathrm{PT}_{\boldsymbol{\mathcal{A}}}]$. We call it the $S$-sorted set of \emph{path term classes}. Let us  note that 
$[\mathrm{PT}_{\boldsymbol{\mathcal{A}}}]$ is a subset of the quotient $\mathrm{T}_{\Sigma^{\boldsymbol{\mathcal{A}}}}(X)/\Theta^{[1]}$, i.e., that $[\mathrm{PT}_{\boldsymbol{\mathcal{A}}}]$ is a subquotient of $\mathrm{T}_{\Sigma^{\boldsymbol{\mathcal{A}}}}(X)$. Actually, we have that
$$
\left[
\mathrm{PT}_{\boldsymbol{\mathcal{A}}}\right] = 
\mathrm{PT}_{\boldsymbol{\mathcal{A}}}/\Theta^{[1]}\!\upharpoonright\!\mathrm{PT}_{\boldsymbol{\mathcal{A}}}.
$$

\index{projection!first-order!$\mathrm{pr}^{\Theta^{[1]}}$}
The projection $\mathrm{pr}^{\Theta^{[1]}}$ from $\mathrm{T}_{\Sigma^{\boldsymbol{\mathcal{A}}}}(X)$ to $\mathrm{T}_{\Sigma^{\boldsymbol{\mathcal{A}}}}(X)/{\Theta^{[1]}}$ birestricts to $\mathrm{PT}_{\boldsymbol{\mathcal{A}}}$ and $[\mathrm{PT}_{\boldsymbol{\mathcal{A}}}]$ and we also let 
$\mathrm{pr}^{\Theta^{[1]}}$ stand for the projection from $\mathrm{PT}_{\boldsymbol{\mathcal{A}}}$ to 
$[\mathrm{PT}_{\boldsymbol{\mathcal{A}}}]$ that, for every sort $s\in S$, maps a  path term $P$ in $\mathrm{PT}_{\boldsymbol{\mathcal{A}},s}$ to $[P]_{s}$, its $\Theta^{[1]}_{s}$-equivalence class, that is
$$
\mathrm{pr}^{\Theta^{[1]}}
\colon
\mathrm{PT}_{\boldsymbol{\mathcal{A}}}
\mor
[\mathrm{PT}_{\boldsymbol{\mathcal{A}}}].
$$
\end{restatable}


We next establish the notation that we will later use for the composition of the embeddings introduced in Definition~\ref{DEta} with the projection to $\Theta^{[1]}$.

\begin{restatable}{definition}{DPTQEta}
\label{DPTQEta}  We will denote by
\begin{enumerate}
\item $\eta^{([1],X)}$ the $S$-sorted mapping from $X$ to 
$[\mathrm{PT}_{\boldsymbol{\mathcal{A}}}]$ given by the composition $\eta^{([1],X)}=\mathrm{pr}^{\Theta^{[1]}}\circ\eta^{(1,X)}$, i.e., for every sort $s\in S$,  $\eta^{([1],X)}$ sends a variable $x\in X_{s}$ to the class $[\eta^{(1,X)}_{s}(x)]_{s}$ in $[\mathrm{PT}_{\boldsymbol{\mathcal{A}}}]_{s}$.
\index{inclusion!first-order!$\eta^{([1],X)}$}
\item $\eta^{([1],\mathcal{A})}$ the $S$-sorted mapping from $\mathcal{A}$ to 
$[\mathrm{PT}_{\boldsymbol{\mathcal{A}}}]$ given by the composition $\eta^{([1],\mathcal{A})}=\mathrm{pr}^{\Theta^{[1]}}\circ\eta^{(1,\mathcal{A})}$, i.e., for every sort $s\in S$,  $\eta^{([1],\mathcal{A})}$ sends a rewrite rule  $\mathfrak{p}\in \mathcal{A}_{s}$ to the class $[\eta^{(1,\mathcal{A})}_{s}(\mathfrak{p})]_{s}$ in $[\mathrm{PT}_{\boldsymbol{\mathcal{A}}}]_{s}$.
\index{inclusion!first-order!$\eta^{([1],\mathcal{A})}$}
\end{enumerate}
The above $S$-sorted mappings are depicted in the diagram of Figure~\ref{FPTQEta}.
\end{restatable}

\begin{figure}
\begin{tikzpicture}
[ACliment/.style={-{To [angle'=45, length=5.75pt, width=4pt, round]}},scale=.8]
\node[] (x) at (0,0) [] {$X$};
\node[] (a) at (0,-1.5) [] {$\mathcal{A}$};
\node[] (T) at (6,-1.5) [] {$
[\mathrm{PT}_{\boldsymbol{\mathcal{A}}}]$};
\draw[ACliment, bend left=10]  (x) to node [above right] {$\eta^{([1],X)}$} (T);
\draw[ACliment]  (a) to node [below] {$\eta^{([1],\mathcal{A})}$} (T);
\end{tikzpicture}
\caption{Quotient embeddings relative to $X$ and $\mathcal{A}$ at layer 1.}\label{FPTQEta}
\end{figure}

We next consider the $S$-sorted mapping from $\mathrm{T}_{\Sigma}(X)$ to 
$[\mathrm{PT}_{\boldsymbol{\mathcal{A}}}]$ obtained from the embedding introduced in Proposition~\ref{PEmb} and the projection $\mathrm{pr}^{\Theta^{[1]}}$.

\begin{restatable}{definition}{DPTQEtas}
\label{DPTQEtas}  We will denote by
\begin{enumerate}
\item  $\eta^{([1],0)\sharp}$ the $S$-sorted mapping from $\mathrm{T}_{\Sigma}(X)$ to 
$[\mathrm{PT}_{\boldsymbol{\mathcal{A}}}]$ that, for every sort $s\in S$, sends a term $P\in \mathrm{T}_{\Sigma}(X)_{s}$ to the  path term class  $[\eta^{(1,0)\sharp}_{s}(P)]_{s}$ in $[\mathrm{PT}_{\boldsymbol{\mathcal{A}}}]_{s}$.
\index{inclusion!first-order!$\eta^{([1],0)\sharp}$}
\end{enumerate}
\end{restatable}

\begin{proposition}\label{PPTQX} The equality $\eta^{([1],X)}=\eta^{([1],0)\sharp}\circ \eta^{(0,X)}$ holds, i.e., the diagram in Figure~\ref{FPTQX} commutes.
\end{proposition}
\begin{proof} 
Let $s$ be a sort in $S$ and $x$ a variable in $X_{s}$. 
The following chain of equalities holds
\allowdisplaybreaks
\begin{align*}
\eta^{([1],0)\sharp}_{s}\left(
\eta^{(0,X)}_{s}\left(
x
\right)
\right)
&=
\left[
\eta^{(1,0)\sharp}_{s}\left(
\eta^{(0,X)}_{s}\left(
x
\right)
\right)
\right]_{s}
\tag{1}
\\&=
\left[
\eta^{(1,X)}_{s}\left(
x
\right)
\right]_{s}
\tag{2}
\\&=
\eta^{([1],X)}_{s}\left(
x
\right).
\tag{3}
\end{align*}

The first  equality applies the mapping $\eta^{([1],0)\sharp}$ according to Definition~\ref{DPTQEtas}; the second equality follows from Proposition~\ref{PEmb}; finally, the last equality recovers the mapping $\eta^{([1],X)}$ introduced in Proposition~\ref{DPTQEta}.
\end{proof}

\begin{figure}
\begin{center}
\begin{tikzpicture}
[ACliment/.style={-{To [angle'=45, length=5.75pt, width=4pt, round]}},scale=0.8]
\node[] (x) at (0,0) [] {$X$};
\node[] (t) at (6,0) [] {$\mathrm{T}_{\Sigma}(X)$};
\node[] (t1) at (6,-2) [] {$[\mathrm{PT}_{\boldsymbol{\mathcal{A}}}]$};

\draw[ACliment]  (x) to node [ above right]
{$\textstyle \eta^{(0,X)}$} (t);
\draw[ACliment, bend right=10]  (x) to node [below left]
{$\textstyle \eta^{([1],X)}$} (t1);

\draw[ACliment] 
(t) to node [right] {$\textstyle \eta^{([1],0)\sharp}$}  (t1);

\end{tikzpicture}
\end{center}
\caption{Quotient embeddings relative to $X$ at layers 0 \& 1.}
\label{FPTQX}
\end{figure}

We next consider the case of the free completion of the $\mathrm{ip}^{(1,X)}$ mapping.

\begin{restatable}{corollary}{CPTQKer}
\label{CPTQKer} 
Let $s$ be a sort in $S$, and $Q$, $P$  path terms in $\mathrm{PT}_{\boldsymbol{\mathcal{A}},s}$ such that $(Q,P)\in\Theta^{[1]}_{s}$. Then 
$
(Q,P)\in\mathrm{Ker}(\mathrm{pr}^{\mathrm{Ker}(\mathrm{CH}^{(1)})}\circ\mathrm{ip}^{(1,X)@}).
$
\end{restatable}
\begin{proof}
It follows from Corollary~\ref{CThetaCong} and Lemma~\ref{LThetaCong}.
\end{proof}

\begin{restatable}{definition}{DPTQIp}
\label{DPTQIp} 
\index{identity!first-order!$\mathrm{ip}^{([1],X)a}$}
We let $\mathrm{ip}^{([1],X)@}$ stand for \[\left(\mathrm{pr}^{\mathrm{Ker}(\mathrm{CH}^{(1)})}\circ\mathrm{ip}^{(1,X)@}\right)^{\mathrm{m}}\circ\mathrm{p}^{ \mathrm{Ker}(\mathrm{pr}^{\mathrm{Ker}(\mathrm{CH}^{(1)})}\circ\mathrm{ip}^{(1,X)@}),\Theta^{[1]}},\]
that is, $\mathrm{ip}^{([1],X)@}$ is the unique many-sorted mapping from $[\mathrm{PT}_{\boldsymbol{\mathcal{A}}}]$ to $[\mathrm{Pth}_{\boldsymbol{\mathcal{A}}}]$
satisfying that
$$
\mathrm{pr}^{\mathrm{Ker}(\mathrm{CH}^{(1)})}\circ\mathrm{ip}^{(1,X)@}
=
\mathrm{ip}^{([1],X)@}
\circ
\mathrm{pr}^{\Theta^{[1]}}.
$$ 
This $S$-sorted mapping is well-defined according to Corollary~\ref{CPTQKer}.
\end{restatable}

We conclude with the Curry-Howard mapping.

\begin{restatable}{definition}{DPTQCH}
\label{DPTQCH} 
\index{Curry-Howard!first-order!$\mathrm{CH}^{[1]}$}
We let $\mathrm{CH}^{[1]}$ stand for the $S$-sorted mapping from $[\mathrm{Pth}_{\boldsymbol{\mathcal{A}}}]$ to $[\mathrm{PT}_{\boldsymbol{\mathcal{A}}}]$ given by the composition $ \mathrm{pr}^{\Theta^{[1]}}\circ\mathrm{CH}^{(1)\mathrm{m}}$, i.e.,
$$
\mathrm{CH}^{[1]}\colon[\mathrm{Pth}_{\boldsymbol{\mathcal{A}}}]\mor
[\mathrm{PT}_{\boldsymbol{\mathcal{A}}}]_{\Theta^{[1]}}.
$$
\end{restatable}

\begin{proposition}\label{PPTQXEq} The following equalities holds
\begin{multicols}{2}
\begin{itemize}
\item[(i)] $\mathrm{ip}^{([1],X)@}\circ\eta^{([1],X)}=\mathrm{ip}^{([1],X)};$
\item[(ii)] $\mathrm{ip}^{([1],X)@}\circ\eta^{([1],0)\sharp}=
\mathrm{ip}^{([1],0)\sharp};
$
\item[(iii)] $\mathrm{CH}^{[1]}\circ\mathrm{ip}^{([1],X)}=\eta^{([1],X)}$;
\item[(iv)] $\mathrm{CH}^{[1]}\circ\mathrm{ip}^{([1],0)\sharp}=\eta^{([1],0)\sharp}$
\end{itemize}
\end{multicols}
The reader is advised to consult the diagram presented in Figure~\ref{FPTQXEq}.
\end{proposition}
\begin{proof}
Regarding the first item, let $s$ be a sort in $S$ and $x$ a variable in $X_{s}$. The following chain of equalities holds
\allowdisplaybreaks
\begin{align*}
\mathrm{ip}^{([1],X)@}_{s}\left(
\eta^{([1],X)}_{s}\left(
x
\right)\right)
&=
\mathrm{ip}^{([1],X)@}_{s}\left(
\left[
\eta^{(1,X)}_{s}\left(
x
\right)
\right]_{s}
\right)
\tag{1}
\\&=
\left[
\mathrm{ip}^{(1,X)@}_{s}\left(
\eta^{(1,X)}_{s}\left(
x
\right)
\right)
\right]_{s}
\tag{2}
\\&=
\left[
\mathrm{ip}^{(1,X)}_{s}\left(
x
\right)
\right]_{s}
\tag{3}
\\&=
\mathrm{ip}^{([1],X)}_{s}\left(
x
\right).
\tag{4}
\end{align*}

The first equality applies the mapping $\eta^{([1],X)}$ introduced in Definition~\ref{DPTQEta}; the second equality applies the mapping $\mathrm{ip}^{([1],X)@}$ introduced in Definition~\ref{DPTQIp}; the third equality follows from Definition~\ref{DIp}; finally, the last equality recovers the description of the mapping $\mathrm{ip}^{([1],X)}$ introduced in Definition~\ref{DCHEch}.

Regarding the second item, let $s$ be a sort in $S$ and $P$ a term in $\mathrm{T}_{\Sigma}(X)_{s}$. The following chain of equalities holds
\allowdisplaybreaks
\begin{align*}
\mathrm{ip}^{([1],X)@}_{s}\left(
\eta^{([1],0)\sharp}_{s}\left(
P
\right)
\right)
&=
\mathrm{ip}^{([1],X)@}_{s}\left(
\left[
\eta^{(1,0)\sharp}_{s}\left(
P
\right)
\right]_{s}
\right)
\tag{1}
\\&=
\left[
\mathrm{ip}^{(1,X)@}_{s}\left(
\eta^{(1,0)\sharp}_{s}\left(
P
\right)
\right)
\right]_{s}
\tag{2}
\\&=
\left[
\mathrm{ip}^{(1,0)\sharp}_{s}\left(
P
\right)
\right]_{s}
\tag{3}
\\&=
\mathrm{ip}^{([1],0)\sharp}_{s}\left(
P
\right).
\tag{4}
\end{align*}

The first equality applies the mapping $\eta^{([1],0)\sharp}$ introduced in Definition~\ref{DPTQEtas}; the second equality applies the mapping $\mathrm{ip}^{([1],X)@}$ introduced in Definition~\ref{DPTQIp}; the third equality follows from Proposition~\ref{PIpUZ}; finally, the last equality recovers the description of the mapping $\mathrm{ip}^{([1],X)}$ introduced in Definition~\ref{DCHEch}.

Regarding the third item, let $s$ be a sort in $S$ and $x$ a variable in $X_{s}$. The following chain of equalities holds
\allowdisplaybreaks
\begin{align*}
\mathrm{CH}^{[1]}_{s}\left(
\mathrm{ip}^{([1],X)}_{s}\left(
x
\right)
\right)
&=
\mathrm{CH}^{[1]}_{s}\left(
\left[
\mathrm{ip}^{(1,X)}_{s}(x)
\right]_{s}
\right)
\tag{1}
\\&=
\left[
\mathrm{CH}^{(1)}_{s}\left(
\mathrm{ip}^{(1,X)}_{s}\left(
x
\right)
\right)
\right]_{s}
\tag{2}
\\&=
\left[
\eta^{(1,X)}_{s}\left(
x
\right)
\right]_{s}
\tag{3}
\\&=
\eta^{([1],X)}_{s}\left(
x
\right).
\tag{4}
\end{align*}

The first equality applies the mapping $\mathrm{ip}^{([1],X)}$ introduced in Definition~\ref{DCHEch}; the second equality applies the mapping $\mathrm{CH}^{[1]}$ introduced in Definition~\ref{DPTQCH}; the third equality follows from Proposition~\ref{PCHId}; finally, the last equality recovers the mapping $\eta^{([1],X)}$ introduced in Definition~\ref{DPTQEtas}.

Regarding the fourth item, let $s$ be a sort in $S$ and $P$ a term in $\mathrm{T}_{\Sigma}(X)_{s}$. The following chain of equalities holds
\allowdisplaybreaks
\begin{align*}
\mathrm{CH}^{[1]}_{s}\left(
\mathrm{ip}^{([1],0)\sharp}_{s}\left(
P
\right)
\right)
&=
\mathrm{CH}^{[1]}_{s}\left(
\left[
\mathrm{ip}^{(1,0)\sharp}_{s}\left(
P
\right)
\right]_{s}
\right)
\tag{1}
\\&=
\left[
\mathrm{CH}^{(1)}_{s}\left(
\mathrm{ip}^{(1,0)\sharp}_{s}\left(
P
\right)
\right)
\right]_{s}
\tag{2}
\\&=
\left[
\eta^{(1,0)\sharp}_{s}\left(
P
\right)
\right]_{s}
\tag{3}
\\&=
\eta^{([1],0)\sharp}_{s}\left(
P
\right).
\tag{4}
\end{align*}

The first equality applies the mapping $\mathrm{ip}^{([1],0)\sharp}$ introduced in Definition~\ref{DCHUZ}; the second equality applies the mapping $\mathrm{CH}^{[1]}$ introduced in Definition~\ref{DPTQCH}; the third equality follows from Proposition~\ref{PCHId}; finally, the last equality recovers the mapping $\eta^{([1],0)\sharp}$ introduced in Definition~\ref{DPTQEtas}.

This completes the proof.
\end{proof}

\begin{figure}
\begin{tikzpicture}
[ACliment/.style={-{To [angle'=45, length=5.75pt, width=4pt, round]}},scale=0.8]
\node[] (x) at (0,0) [] {$X$};
\node[] (t) at (6,0) [] {$\mathrm{T}_{\Sigma}(X)$};
\node[] (pth1) at (6,-2) [] {$[\mathrm{Pth}_{\boldsymbol{\mathcal{A}}}]$};
\node[] (pt1) at (6,-4) [] {$[\mathrm{PT}_{\boldsymbol{\mathcal{A}}}]$};

\draw[ACliment]  (x) to node [above right]
{$\textstyle \eta^{(0,X)}$} (t);
\draw[ACliment, bend right=10]  (x) to node [pos=.5, fill=white]
{$\textstyle \mathrm{ip}^{([1],X)}$} (pth1.west);
\draw[ACliment, bend right=20]  (x) to node  [below left]
{$\textstyle \eta^{([1],X)}$} (pt1.west);

\draw[ACliment] 
($(t)+(0,-.35)$) to node [right] {$\textstyle \mathrm{ip}^{([1],0)\sharp}$}  ($(pth1)+(0,.35)$);

\draw[ACliment] 
($(pth1)+(-.30,-.35)$) to node [left] {$\textstyle \mathrm{CH}^{[1]}$}  ($(pt1)+(-.30,.35)$);
\draw[ACliment] 
($(pt1)+(.30,+.35)$) to node [right] {$\textstyle \mathrm{ip}^{([1],X)@}$}  ($(pth1)+(.30,-.35)$);

\draw[ACliment, rounded corners]
(t.east) -- 
($(t)+(2.5,0)$) -- node [right] {$\textstyle \eta^{([1],0)\sharp}$} 
($(pt1)+(2.5,0)$) -- (pt1.east)
;

\end{tikzpicture}
\caption{Quotient embeddings relative to $X$ at layers $0$ \& $1$.}
\label{FPTQXEq}
\end{figure}

\begin{proposition}\label{PPTQAEq} The following equalities are fulfilled
\begin{multicols}{2}
\begin{itemize}
\item[(i)] $\mathrm{ip}^{([1],X)@}\circ\eta^{([1],\mathcal{A})}=\mathrm{ech}^{([1],\mathcal{A})};$
\item[(ii)] $\mathrm{CH}^{[1]}\circ\mathrm{ech}^{([1],\mathcal{A})}=\eta^{([1],\mathcal{A})}$.
\end{itemize}
\end{multicols}
The reader is advised to consult the diagram presented in Figure~\ref{FPTQAEq}.
\end{proposition}

\begin{figure}
\begin{tikzpicture}
[ACliment/.style={-{To [angle'=45, length=5.75pt, width=4pt, round]}},scale=0.8]
\node[] (a) at (0,-2) [] {$\mathcal{A}$};
\node[] (pth1) at (6,-2) [] {$[\mathrm{Pth}_{\boldsymbol{\mathcal{A}}}]$};
\node[] (pt1) at (6,-4) [] {$[\mathrm{PT}_{\boldsymbol{\mathcal{A}}}]$};

\draw[ACliment]  (a) to node [above]
{$\textstyle \mathrm{ech}^{([1],\mathcal{A})}$} (pth1);
\draw[ACliment, bend right=10]  (a) to node  [below left]
{$\textstyle \eta^{([1],\mathcal{A})}$} (pt1.west);

\draw[ACliment] 
($(pth1)+(-.3,-.35)$) to node [left] {$\textstyle \mathrm{CH}^{[1]}$}  ($(pt1)+(-.3,.35)$);
\draw[ACliment] 
($(pt1)+(.3,+.35)$) to node [right] {$\textstyle \mathrm{ip}^{([1],X)@}$}  ($(pth1)+(.3,-.35)$);
\end{tikzpicture}
\caption{Quotient embeddings relative to $\mathcal{A}$ at layer $1$.}
\label{FPTQAEq}
\end{figure}

\section{
\texorpdfstring
{A structure of partial $\Sigma^{\boldsymbol{\mathcal{A}}}$-algebra on $[\mathrm{PT}_{\boldsymbol{\mathcal{A}}}]$}
{A partial algebra on the quotient of path terms}
}

In this section we show that the $S$-sorted set of path term classes is equipped with a structure of partial $\Sigma^{\boldsymbol{\mathcal{A}}}$-algebra so that the resulting algebra is a $\Sigma^{\boldsymbol{\mathcal{A}}}$-subalgebra of $\mathbf{T}_{\Sigma^{\boldsymbol{\mathcal{A}}}}(X)/{\Theta^{[1]}}$. 

\begin{restatable}{proposition}{PPTQCatAlg}
\label{PPTQCatAlg} 
\index{path terms!first-order!$[\mathbf{PT}_{\boldsymbol{\mathcal{A}}}]$}
The $S$-sorted set $[\mathrm{PT}_{\boldsymbol{\mathcal{A}}}]$ is equipped, in a natural way, with a structure of many-sorted partial $\Sigma^{\boldsymbol{\mathcal{A}}}$-algebra, which is a 
$\Sigma^{\boldsymbol{\mathcal{A}}}$-subalgebra of $\mathbf{T}_{\Sigma^{\boldsymbol{\mathcal{A}}}}(X)/{\Theta^{[1]}}$.
\end{restatable}
\begin{proof}
Let us denote by $[\mathbf{PT}_{\boldsymbol{\mathcal{A}}}]$ the partial $\Sigma^{\boldsymbol{\mathcal{A}}}$-algebra defined as follows:

\textsf{(i)} The underlying $S$-sorted set of $[\mathbf{PT}_{\boldsymbol{\mathcal{A}}}]$ is $[\mathrm{PT}_{\boldsymbol{\mathcal{A}}}]$.

\textsf{(ii)} For every $(\mathbf{s},s)\in S^{\star}\times S$ and every operation symbol $\sigma\in\Sigma_{\mathbf{s},s}$, the operation $\sigma^{[\mathbf{PT}_{\boldsymbol{\mathcal{A}}}]}$ is equal to $\sigma^{\mathbf{T}_{\Sigma^{\boldsymbol{\mathcal{A}}}}(X)/{\Theta^{[1]}}}$, i.e., to the interpretation of $\sigma$ in $\mathbf{T}_{\Sigma^{\boldsymbol{\mathcal{A}}}}(X)/{\Theta^{[1]}}$. 

Let us note that, if $([P_{j}]_{s_{j}})_{j\in\bb{\mathbf{s}}}$ is a family of path term classes in $[\mathrm{PT}_{\boldsymbol{\mathcal{A}}}]_{{\mathbf{s}}}$, then the term class
$$
\sigma^{[\mathbf{PT}_{\boldsymbol{\mathcal{A}}}]}
\left(\left(
\left[P_{j}\right]_{s_{j}}
\right)_{j\in\bb{\mathbf{s}}}\right)
=
\left[
\sigma^{\mathbf{T}_{\Sigma^{\boldsymbol{\mathcal{A}}}}(X)}
\left(\left(P_{j}\right)_{j\in\bb{\mathbf{s}}}\right)
\right]_{s}
$$
is a path term class by Proposition~\ref{PPT}. That is, the operation $\sigma^{[\mathbf{PT}_{\boldsymbol{\mathcal{A}}}]}$ is well-defined.

\textsf{(iii)}  For every $s\in S$ and every rewrite rule $\mathfrak{p}\in \mathcal{A}_{s}$, the constant $\mathfrak{p}^{[\mathbf{PT}_{\boldsymbol{\mathcal{A}}}]}$ is equal to $\mathfrak{p}^{\mathbf{T}_{\Sigma^{\boldsymbol{\mathcal{A}}}}(X)/{\Theta^{[1]}}}$, i.e., to the interpretation of $\mathfrak{p}$ in $\mathbf{T}_{\Sigma^{\boldsymbol{\mathcal{A}}}}(X)/{\Theta^{[1]}}$. 

Let us note that the term class
$$
\mathfrak{p}^{[\mathbf{PT}_{\boldsymbol{\mathcal{A}}}]}
=
\left[
\mathfrak{p}^{\mathbf{T}_{\Sigma^{\boldsymbol{\mathcal{A}}}}(X)}
\right]_{s}
$$
is a path term class by Proposition~\ref{PPT}. That is, the constant $
\mathfrak{p}^{[\mathbf{PT}_{\boldsymbol{\mathcal{A}}}]}$ is well-defined.

\textsf{(iv)}  For every $s\in S$, the operation $\mathrm{sc}_{s}^{0[\mathbf{PT}_{\boldsymbol{\mathcal{A}}}]}$ is equal to $\mathrm{sc}_{s}^{0\mathbf{T}_{\Sigma^{\boldsymbol{\mathcal{A}}}}(X)/{\Theta^{[1]}}}$, i.e., to the interpretation of $\mathrm{sc}_{s}^{0}$ in $\mathbf{T}_{\Sigma^{\boldsymbol{\mathcal{A}}}}(X)/{\Theta^{[1]}}$. 

Let us note that, if $[P]_{s}$ is a  path term class in $[\mathrm{PT}_{\boldsymbol{\mathcal{A}}}]_{s}$, then the term class
$$
\mathrm{sc}_{s}^{0[\mathbf{PT}_{\boldsymbol{\mathcal{A}}}]}
\left(
\left[P\right]_{s}
\right)
=
\left[
\mathrm{sc}_{s}^{0\mathbf{T}_{\Sigma^{\boldsymbol{\mathcal{A}}}}(X)}
(P)
\right]_{s}
$$
is a path term class by Proposition~\ref{PPT}. That is, the operation $\mathrm{sc}_{s}^{0[\mathbf{PT}_{\boldsymbol{\mathcal{A}}}]}$ is well-defined.

\textsf{(v)}  For every $s\in S$, the operation $\mathrm{tg}_{s}^{0[\mathbf{PT}_{\boldsymbol{\mathcal{A}}}]}$ is equal to $\mathrm{tg}_{s}^{0\mathbf{T}_{\Sigma^{\boldsymbol{\mathcal{A}}}}(X)/{\Theta^{[1]}}}$, i.e., to the interpretation of $\mathrm{tg}_{s}^{0}$ in $\mathbf{T}_{\Sigma^{\boldsymbol{\mathcal{A}}}}(X)/{\Theta^{[1]}}$.  The operation $\mathrm{tg}_{s}^{0[\mathbf{PT}_{\boldsymbol{\mathcal{A}}}]}$ is well-defined by a similar argument to that of $\mathrm{sc}_{s}^{0[\mathbf{PT}_{\boldsymbol{\mathcal{A}}}]}$.

\textsf{(vi)}  For every $s\in S$, the partial binary operation $\circ_{s}^{0[\mathbf{PT}_{\boldsymbol{\mathcal{A}}}]}$ is defined on path term classes $[Q]_{s}$, $[P]_{s}$ in $\mathrm{PT}_{\boldsymbol{\mathcal{A}},s}$ if, and only if, the following equation holds
$$
\mathrm{sc}_{s}^{(0,1)}\left(\mathrm{ip}^{(1,X)@}_{s}\left(Q\right)\right)=
\mathrm{tg}_{s}^{(0,1)}\left(\mathrm{ip}^{(1,X)@}_{s}\left(P\right)\right).
$$

Note that the last equation does not depend on the representative in the $\Theta^{[1]}_{s}$-class by Corollary~\ref{CPTScTg}.

For this case, $\circ_{s}^{0[\mathbf{PT}_{\boldsymbol{\mathcal{A}}}]}$ is equal to $\circ_{s}^{\mathbf{T}_{\Sigma^{\boldsymbol{\mathcal{A}}}}(X)/{\Theta^{[1]}}}$, i.e., to the interpretation of $\circ^{0}_{s}$ in $\mathbf{T}_{\Sigma^{\boldsymbol{\mathcal{A}}}}(X)/{\Theta^{[1]}}$. 

Let us note that, if $[Q]_{s}$, $[P]_{s}$ are  path term classes in $[\mathrm{PT}_{\boldsymbol{\mathcal{A}}}]_{s}$ with
$
\mathrm{sc}_{s}^{(0,1)}(\mathrm{ip}^{(1,X)@}_{s}(Q))=
\mathrm{tg}_{s}^{(0,1)}(\mathrm{ip}^{(1,X)@}_{s}(P)),
$ 
then the term class
$$
[Q]_{s}
\circ_{s}^{0[\mathbf{PT}_{\boldsymbol{\mathcal{A}}}]}
[P]_{s}
=
\left[
Q
\circ_{s}^{0\mathbf{T}_{\Sigma^{\boldsymbol{\mathcal{A}}}}(X)}
P
\right]_{s}
$$
is a path term class by Proposition~\ref{PPTCatAlg}. That is, the operation $\circ_{s}^{0[\mathbf{PT}_{\boldsymbol{\mathcal{A}}}]}$ is well-defined.

This completes the definition of $[\mathbf{PT}_{\boldsymbol{\mathcal{A}}}]$, the partial many-sorted $\Sigma^{\boldsymbol{\mathcal{A}}}$-algebra of path term classes. Thus defined it is obvious that 
$[\mathbf{PT}_{\boldsymbol{\mathcal{A}}}]$ is a $\Sigma^{\boldsymbol{\mathcal{A}}}$-subalgebra of $\mathbf{T}_{\Sigma^{\boldsymbol{\mathcal{A}}}}(X)/{\Theta^{[1]}}$.
\end{proof}

We next state that for the many-sorted partial $\Sigma^{\boldsymbol{\mathcal{A}}}$-algebras $\mathbf{PT}_{\boldsymbol{\mathcal{A}}}$, of path terms, and $[\mathbf{PT}_{\boldsymbol{\mathcal{A}}}]$, of  path term classes, introduced in Proposition~\ref{PPTCatAlg} and Proposition~\ref{PPTQCatAlg}, respectively, the canonical projection $\mathrm{pr}^{\Theta^{[1]}}$ determines a surjective 
$\Sigma^{\boldsymbol{\mathcal{A}}}$-homomorphism. 

\begin{restatable}{corollary}{CPTQPr}
\label{CPTQPr} The mapping $\mathrm{pr}^{\Theta^{[1]}}$ is a surjective $\Sigma^{\boldsymbol{\mathcal{A}}}$-homomorphism
$$
\mathrm{pr}^{\Theta^{[1]}}\colon
\mathbf{PT}_{\boldsymbol{\mathcal{A}}}
\mor
[\mathbf{PT}_{\boldsymbol{\mathcal{A}}}].
$$
\end{restatable}

We finally prove a lemma that will be of use later. It states that the operations coming from $\Sigma$ are injective in $[\mathbf{PT}_{\boldsymbol{\mathcal{A}}}]$.

\begin{lemma}\label{LPTQSigma} 
Let $\mathbf{s}$ be a word in $S^{\star}$, $s$ a sort in $S$, $\sigma$ an operation symbol in $\Sigma_{\mathbf{s},s}$, and $(P_{j})_{j\in\bb{\mathbf{s}}}$, $(Q_{j})_{j\in\bb{\mathbf{s}}}$ two families of path terms in $\mathrm{PT}_{\boldsymbol{\mathcal{A}},\mathbf{s}}$. If the following equality is fulfilled
$$
\sigma^{[\mathbf{PT}_{\boldsymbol{\mathcal{A}}}]}\left(\left(\left[P_{j}
\right]_{s_{j}}\right)_{j\in\bb{\mathbf{s}}}\right)
=
\sigma^{[\mathbf{PT}_{\boldsymbol{\mathcal{A}}}]}\left(\left(\left[Q_{j}
\right]_{s_{j}}\right)_{j\in\bb{\mathbf{s}}}\right),
$$
then, for every $j\in\bb{\mathbf{s}}$, we have that 
$
[P_{j}]_{s_{j}}
=
[Q_{j}]_{s_{j}}.
$
\end{lemma}
\begin{proof}
By assumption we have that 
$$
\left[
\sigma^{\mathbf{PT}_{\boldsymbol{\mathcal{A}}}}
\left(\left(
P_{j}\right)_{j\in\bb{\mathbf{s}}}\right)
\right]_{s}
=
\left[
\sigma^{\mathbf{PT}_{\boldsymbol{\mathcal{A}}}}
\left(\left(
Q_{j}\right)_{j\in\bb{\mathbf{s}}}\right)
\right]_{s}.
$$
Thus, by Lemma~\ref{LThetaCong}, we have that 
\begin{multline*}
\mathrm{CH}^{(1)}_{s}\left(
\mathrm{ip}^{(1,X)@}_{s}\left(
\sigma^{\mathbf{PT}_{\boldsymbol{\mathcal{A}}}}
\left(\left(
P_{j}\right)_{j\in\bb{\mathbf{s}}}\right)
\right)\right)
\\=
\mathrm{CH}^{(1)}_{s}\left(\mathrm{ip}^{(1,X)@}_{s}\left(
\sigma^{\mathbf{PT}_{\boldsymbol{\mathcal{A}}}}
\left(\left(
Q_{j}\right)_{j\in\bb{\mathbf{s}}}\right)
\right)\right).
\tag{3}
\end{multline*}

The following chain of equalities holds
\begin{flushleft}
$
\mathrm{CH}^{(1)}_{s}\left(
\mathrm{ip}^{(1,X)@}_{s}\left(
\sigma^{\mathbf{PT}_{\boldsymbol{\mathcal{A}}}}
\left(\left(
P_{j}\right)_{j\in\bb{\mathbf{s}}}\right)
\right)\right)$
\allowdisplaybreaks
\begin{align*}
\qquad&=
\mathrm{CH}^{(1)}_{s}\left(
\sigma^{\mathbf{F}_{\Sigma^{\boldsymbol{\mathcal{A}}}}
\left(\mathbf{Pth}_{\boldsymbol{\mathcal{A}}}\right)}
\left(\left(
\mathrm{ip}^{(1,X)@}_{s_{j}}\left(
P_{j}
\right)\right)_{j\in\bb{\mathbf{s}}}\right)\right)
\tag{1}
\\&=
\mathrm{CH}^{(1)}_{s}\left(
\sigma^{\mathbf{Pth}_{\boldsymbol{\mathcal{A}}}}\left(\left(
\mathrm{ip}^{(1,X)@}_{s_{j}}\left(
P_{j}
\right)\right)_{j\in\bb{\mathbf{s}}}
\right)\right)
\tag{2}
\\&=
\sigma^{\mathbf{T}_{\Sigma^{\boldsymbol{\mathcal{A}}}}(X)}
\left(\left(
\mathrm{CH}^{(1)}_{s_{j}}\left(
\mathrm{ip}^{(1,X)@}_{s_{j}}\left(
P_{j}
\right)\right)\right)_{j\in\bb{\mathbf{s}}}\right).
\tag{3}
\end{align*}
\end{flushleft}

The first equality follows from the fact that, by Definition~\ref{DIp}, $\mathrm{ip}^{(1,X)@}$ is a $\Sigma^{\boldsymbol{\mathcal{A}}}$-homomorphism; the second equality follows from the fact that, by assumption, for every $j\in\bb{\mathbf{s}}$, $P_{j}$ is a path term. By Proposition~\ref{PPT}, for every $j\in\bb{\mathbf{s}}$, $\mathrm{ip}^{(1,X)@}_{s_{j}}(P_{j})$ is a path in $\mathbf{Pth}_{\boldsymbol{\mathcal{A}},s_{j}}$. Hence, the interpretation of the operation symbol $\sigma$ in the free many-sorted $\Sigma^{\boldsymbol{\mathcal{A}}}$-algebra $\mathbf{F}_{\Sigma^{\boldsymbol{\mathcal{A}}}}(\mathbf{Pth}_{\boldsymbol{\mathcal{A}}})$ applied to the family $(\mathrm{ip}^{(1,X)@}_{s_{j}}(P_{j}))_{j\in\bb{\mathbf{s}}}$ is given by the interpretation of the operation symbol $\sigma$ in the partial many-sorted $\Sigma^{\boldsymbol{\mathcal{A}}}$-algebra $\mathbf{Pth}_{\boldsymbol{\mathcal{A}}}$ on this same family; finally, the last equality follows from the fact that, according to Proposition~\ref{PCHHom}, the Curry-Howard mapping is a $\Sigma$-homomorphism.

By a similar argument we also have that 
\begin{multline*}
\mathrm{CH}^{(1)}_{s}\left(
\mathrm{ip}^{(1,X)@}_{s}\left(
\sigma^{\mathbf{PT}_{\boldsymbol{\mathcal{A}}}}
\left(\left(
Q_{j}\right)_{j\in\bb{\mathbf{s}}}
\right)\right)\right)
\\=
\sigma^{\mathbf{T}_{\Sigma^{\boldsymbol{\mathcal{A}}}}(X)}
\left(\left(
\mathrm{CH}^{(1)}_{s_{j}}\left(
\mathrm{ip}^{(1,X)@}_{s_{j}}\left(
Q_{j}
\right)\right)\right)_{j\in\bb{\mathbf{s}}}\right).
\end{multline*}

Therefore, for every $j\in\bb{\mathbf{s}}$, we have that 
$$
\mathrm{CH}^{(1)}_{s_{j}}\left(
\mathrm{ip}^{(1,X)@}_{s_{j}}\left(
P_{j}
\right)\right)
=
\mathrm{CH}^{(1)}_{s_{j}}\left(
\mathrm{ip}^{(1,X)@}_{s_{j}}\left(
Q_{j}
\right)\right).
$$

According to Lemma~\ref{LWCong}, we have that, for every $j\in\bb{\mathbf{s}}$,
\begin{align*}
\left(P_{j},
\mathrm{CH}^{(1)}_{s_{j}}\left(
\mathrm{ip}^{(1,X)@}_{s_{j}}\left(
P_{j}
\right)\right)\right)\in\Theta^{[1]}_{s_{j}},
&&
\left(Q_{j},
\mathrm{CH}^{(1)}_{s_{j}}\left(
\mathrm{ip}^{(1,X)@}_{s_{j}}\left(
Q_{j}
\right)\right)\right)\in\Theta^{[1]}_{s_{j}}.
\end{align*}

Since $\Theta^{[1]}$ is an equivalence relation, we conclude by transitivity that, for every $j\in\bb{\mathbf{s}}$, 
$
[P_{j}]_{s_{j}}
=
[Q_{j}]_{s_{j}}.
$

This concludes the proof.
\end{proof}


\section{
\texorpdfstring
{A structure of $S$-sorted categorial $\Sigma$-algebra on $[\mathrm{PT}_{\boldsymbol{\mathcal{A}}}]$}
{A categorial algebra on the quotient of path terms}
} 
In what follows we will prove that the $S$-sorted set $[\mathrm{PT}_{\boldsymbol{\mathcal{A}}}]$ is equipped with a structure of $S$-sorted categorial $\Sigma$-algebra. To do so we begin by showing that the partial $\Sigma^{\boldsymbol{\mathcal{A}}}$-algebra $[\mathrm{PT}_{\boldsymbol{\mathcal{A}}}]$ satisfies the defining equations of the concept of $S$-sorted category.


\begin{proposition}\label{PPTQVarA2} Let $s$ be a sort in $S$ and $[P]_{s}$ a path term class in $[\mathrm{PT}_{\boldsymbol{\mathcal{A}}}]_{s}$, then the following equalities holds
\begin{align*}
\mathrm{sc}^{0
[\mathbf{PT}_{\boldsymbol{\mathcal{A}}}]
}_{s}\left(
\mathrm{sc}^{0
[\mathbf{PT}_{\boldsymbol{\mathcal{A}}}]
}_{s}\left(
[P]_{s}
\right)\right)
&=
\mathrm{sc}^{0
[\mathbf{PT}_{\boldsymbol{\mathcal{A}}}]
}_{s}\left(
[P]_{s}
\right);
\\
\mathrm{sc}^{0
[\mathbf{PT}_{\boldsymbol{\mathcal{A}}}]
}_{s}\left(
\mathrm{tg}^{0
[\mathbf{PT}_{\boldsymbol{\mathcal{A}}}]
}_{s}\left(
[P]_{s}
\right)\right)
&=
\mathrm{tg}^{0
[\mathbf{PT}_{\boldsymbol{\mathcal{A}}}]
}_{s}\left(
[P]_{s}
\right);
\\
\mathrm{tg}^{0
[\mathbf{PT}_{\boldsymbol{\mathcal{A}}}]
}_{s}\left(
\mathrm{sc}^{0
[\mathbf{PT}_{\boldsymbol{\mathcal{A}}}]
}_{s}\left(
[P]_{s}
\right)\right)
&=
\mathrm{sc}^{0
[\mathbf{PT}_{\boldsymbol{\mathcal{A}}}]
}_{s}\left(
[P]_{s}
\right);
\\
\mathrm{tg}^{0
[\mathbf{PT}_{\boldsymbol{\mathcal{A}}}]
}_{s}\left(
\mathrm{tg}^{0
[\mathbf{PT}_{\boldsymbol{\mathcal{A}}}]
}_{s}\left(
[P]_{s}
\right)\right)
&=
\mathrm{tg}^{0
[\mathbf{PT}_{\boldsymbol{\mathcal{A}}}]
}_{s}\left(
[P]_{s}
\right).
\end{align*}
\end{proposition}
\begin{proof}

Regarding the first equation, the following chain of equalities holds
\allowdisplaybreaks
\begin{align*}
\mathrm{sc}^{0
[\mathbf{PT}_{\boldsymbol{\mathcal{A}}}]
}_{s}\left(
\mathrm{sc}^{0
[\mathbf{PT}_{\boldsymbol{\mathcal{A}}}]
}_{s}\left(
[P]_{s}
\right)\right)&=
\left[
\mathrm{sc}^{0\mathbf{PT}_{\boldsymbol{\mathcal{A}}}}_{s}\left(
\mathrm{sc}^{0\mathbf{PT}_{\boldsymbol{\mathcal{A}}}}_{s}\left(
P
\right)
\right)
\right]_{s}
\tag{1}
\\&=
\left[
\mathrm{CH}^{(1)}_{s}\left(
\mathrm{ip}^{(1,X)@}_{s}\left(
\mathrm{sc}^{0\mathbf{PT}_{\boldsymbol{\mathcal{A}}}}_{s}\left(
\mathrm{sc}^{0\mathbf{PT}_{\boldsymbol{\mathcal{A}}}}_{s}\left(
P
\right)
\right)
\right)
\right)
\right]_{s}
\tag{2}
\\&=
\left[
\mathrm{CH}^{(1)}_{s}\left(
\mathrm{sc}^{0\mathbf{Pth}_{\boldsymbol{\mathcal{A}}}}_{s}\left(
\mathrm{sc}^{0\mathbf{Pth}_{\boldsymbol{\mathcal{A}}}}_{s}\left(
\mathrm{ip}^{(1,X)@}_{s}\left(
P
\right)
\right)
\right)
\right)
\right]_{s}
\tag{3}
\\&=
\left[
\mathrm{CH}^{(1)}_{s}\left(
\mathrm{sc}^{0\mathbf{Pth}_{\boldsymbol{\mathcal{A}}}}_{s}\left(
\mathrm{ip}^{(1,X)@}_{s}\left(
P
\right)
\right)
\right)
\right]_{s}
\tag{4}
\\&=
\left[
\mathrm{CH}^{(1)}_{s}\left(
\mathrm{ip}^{(1,X)@}_{s}\left(
\mathrm{sc}^{0\mathbf{Pth}_{\boldsymbol{\mathcal{A}}}}_{s}\left(
P
\right)
\right)
\right)
\right]_{s}
\tag{5}
\\&=
\left[
\mathrm{sc}^{0\mathbf{Pth}_{\boldsymbol{\mathcal{A}}}}_{s}\left(
P
\right)
\right]_{s}
\tag{6}
\\&=
\mathrm{sc}^{0
[\mathbf{PT}_{\boldsymbol{\mathcal{A}}}]
}_{s}\left(
[P]_{s}
\right).
\tag{7}
\end{align*}

The first equality unravels the interpretation of the $0$-source operation in the many-sorted partial $\Sigma^{\boldsymbol{\mathcal{A}}}$-algebra $[\mathbf{PT}_{\boldsymbol{\mathcal{A}}}]$ according to Proposition~\ref{PPTQCatAlg}; the second equality follows from Lemma~\ref{LWCong}; the third equality follows from the fact that $\mathrm{ip}^{(1,X)@}$ is a $\Sigma^{\boldsymbol{\mathcal{A}}}$-homomorphism according to Definition~\ref{DIp}; the fourth equality follows from Proposition~\ref{PCHVarA2}; the fifth equality follows from the fact that $\mathrm{ip}^{(1,X)@}$ is a $\Sigma^{\boldsymbol{\mathcal{A}}}$-homomorphism according to Definition~\ref{DIp}; the sixth equality follows from Lemma~\ref{LWCong}; finally, the last equality recovers the interpretation of the  $0$-source operation in the many-sorted partial $\Sigma^{\boldsymbol{\mathcal{A}}}$-algebra $[\mathbf{PT}_{\boldsymbol{\mathcal{A}}}]$ according to Proposition~\ref{PPTQCatAlg}.

The other equalities holds by a similar argument.

This completes the proof.
\end{proof}

\begin{proposition}\label{PPTQVarA3} Let $s$ be a sort in $S$ and $[P]_{s}$, $[Q]_{s}$ path term classes in $[\mathrm{PT}_{\boldsymbol{\mathcal{A}}}]_{s}$, then the following statements are equivalent
\begin{enumerate}
\item[(i)] $[Q]_{s}
\circ^{0
[\mathbf{PT}_{\boldsymbol{\mathcal{A}}}]}_{s}
[P]_{s}
$ is defined;
\item[(ii)] $\mathrm{sc}^{0
[\mathbf{PT}_{\boldsymbol{\mathcal{A}}}]}_{s}(
[Q]_{s}
)
=
\mathrm{tg}^{0
[\mathbf{PT}_{\boldsymbol{\mathcal{A}}}]}_{s}(
[P]_{s}
)
$.
\end{enumerate}
\end{proposition}
\begin{proof}
The following chain of equivalences holds
\begin{flushleft}
$[Q]_{s}
\circ^{0
[\mathbf{PT}_{\boldsymbol{\mathcal{A}}}]}_{s}
[P]_{s}$ is defined
\allowdisplaybreaks
\begin{align*}
&\Leftrightarrow
\mathrm{sc}^{(0,1)}_{s}\left(
\mathrm{ip}^{(1,X)@}_{s}\left(
Q
\right)\right)
=
\mathrm{tg}^{(0,1)}_{s}\left(
\mathrm{ip}^{(1,X)@}_{s}\left(
P
\right)\right)
\tag{1}
\\&\Leftrightarrow
\mathrm{sc}^{0
[\mathbf{Pth}_{\boldsymbol{\mathcal{A}}}]
}_{s}\left(
\left[
\mathrm{ip}^{(1,X)@}_{s}\left(
Q
\right)
\right]_{s}
\right)
=
\mathrm{tg}^{0
[\mathbf{Pth}_{\boldsymbol{\mathcal{A}}}]
}_{s}\left(
\left[
\mathrm{ip}^{(1,X)@}_{s}\left(
P
\right)
\right]_{s}
\right)
\tag{2}
\\&\Leftrightarrow
\left[
\mathrm{sc}^{0\mathbf{Pth}_{\boldsymbol{\mathcal{A}}}}_{s}\left(
\mathrm{ip}^{(1,X)@}_{s}\left(
Q
\right)
\right)
\right]_{s}
=
\left[
\mathrm{tg}^{0\mathbf{Pth}_{\boldsymbol{\mathcal{A}}}}_{s}\left(
\mathrm{ip}^{(1,X)@}_{s}\left(
P
\right)
\right)
\right]_{s}
\tag{3}
\\&\Leftrightarrow
\left[
\mathrm{ip}^{(1,X)@}_{s}\left(
\mathrm{sc}^{0\mathbf{PT}_{\boldsymbol{\mathcal{A}}}}_{s}\left(
Q
\right)
\right)
\right]_{s}
=
\left[
\mathrm{ip}^{(1,X)@}_{s}\left(
\mathrm{tg}^{0\mathbf{PT}_{\boldsymbol{\mathcal{A}}}}_{s}\left(
P
\right)
\right)
\right]_{s}
\tag{4}
\\&\Leftrightarrow
\mathrm{CH}^{(1)}_{s}\left(
\mathrm{ip}^{(1,X)@}_{s}\left(
\mathrm{sc}^{0\mathbf{PT}_{\boldsymbol{\mathcal{A}}}}_{s}\left(
Q
\right)\right)\right)
=
\mathrm{CH}^{(1)}_{s}\left(
\mathrm{ip}^{(1,X)@}_{s}\left(
\mathrm{tg}^{0\mathbf{PT}_{\boldsymbol{\mathcal{A}}}}_{s}\left(
P
\right)\right)\right)
\tag{5}
\\&\Leftrightarrow
\left[
\mathrm{sc}^{0\mathbf{PT}_{\boldsymbol{\mathcal{A}}}}_{s}\left(
Q
\right)
\right]_{s}
=
\left[
\mathrm{tg}^{0\mathbf{PT}_{\boldsymbol{\mathcal{A}}}}_{s}\left(
P
\right)
\right]_{s}
\tag{6}
\\&\Leftrightarrow
\mathrm{sc}^{0
[\mathbf{PT}_{\boldsymbol{\mathcal{A}}}]}_{s}\left(
\left[
Q
\right]_{s}
\right)
=
\mathrm{tg}^{0
[\mathbf{PT}_{\boldsymbol{\mathcal{A}}}]}_{s}\left(
\left[
P
\right]_{s}
\right).
\tag{7}
\end{align*}
\end{flushleft}

The first equivalence follows from the interpretation of the $0$-composition operation in $[\mathbf{PT}_{\boldsymbol{\mathcal{A}}}]$ according to Proposition~\ref{PPTQCatAlg}; the second equivalence follows from Proposition~\ref{PCHVarA3}; the third equivalence follows from the interpretation of the $0$-source and $0$-target operation in $[\mathbf{Pth}_{\boldsymbol{\mathcal{A}}}]$ according to Proposition~\ref{PCHCatAlg}; the fourth equivalence follows from the the fact that $\mathrm{ip}^{(1,X)@}$ is a $\Sigma^{\boldsymbol{\mathcal{A}}}$-homomorphism according to Definition~\ref{DIp}; the fifth equivalence simply unravels the description of a class under the kernel of the Curry-Howard mapping; the sixth equivalence follows, from left to right by Lemma~\ref{LWCong} and from right to left by Lemma~\ref{LThetaCong}; finally, the last equivalence recovers interpretation of the $0$-source and $0$-target operation in $[\mathbf{PT}_{\boldsymbol{\mathcal{A}}}]$ according to Proposition~\ref{PPTQCatAlg}.

This completes the proof.
\end{proof}

\begin{proposition}\label{PPTQVarA4} Let $s$ be a sort in $S$ and $[P]_{s}$, $[Q]_{s}$ path term classes in $[\mathrm{PT}_{\boldsymbol{\mathcal{A}}}]_{s}$.
If $\mathrm{sc}^{0
[\mathbf{PT}_{\boldsymbol{\mathcal{A}}}]
}_{s}([Q]_{s})=
\mathrm{tg}^{0
[\mathbf{PT}_{\boldsymbol{\mathcal{A}}}]
}_{s}([P]_{s})
$, then the following equalities hold
\allowdisplaybreaks
\begin{align*}
\mathrm{sc}^{0
[\mathbf{PT}_{\boldsymbol{\mathcal{A}}}]
}_{s}\left(
[Q]_{s}
\circ^{0
[\mathbf{PT}_{\boldsymbol{\mathcal{A}}}]
}_{s}
[P]_{s}
\right)
&=
\mathrm{sc}^{0
[\mathbf{PT}_{\boldsymbol{\mathcal{A}}}]
}_{s}\left(
[P]_{s}
\right);
\\
\mathrm{tg}^{0
[\mathbf{PT}_{\boldsymbol{\mathcal{A}}}]
}_{s}\left(
[Q]_{s}
\circ^{0
[\mathbf{PT}_{\boldsymbol{\mathcal{A}}}]
}_{s}
[P]_{s}
\right)
&=
\mathrm{tg}^{0
[\mathbf{PT}_{\boldsymbol{\mathcal{A}}}]
}_{s}\left(
[Q]_{s}
\right).
\end{align*}
\end{proposition}
\begin{proof}
The following chain of equalities holds
\begin{flushleft}
$\mathrm{sc}^{0
[\mathbf{PT}_{\boldsymbol{\mathcal{A}}}]
}_{s}\left(
[Q]_{s}
\circ^{0
[\mathbf{PT}_{\boldsymbol{\mathcal{A}}}]
}_{s}
[P]_{s}
\right)$
\allowdisplaybreaks
\begin{align*}
\qquad&=
\left[
\mathrm{sc}^{0\mathbf{PT}_{\boldsymbol{\mathcal{A}}}}_{s}\left(
Q
\circ^{0\mathbf{PT}_{\boldsymbol{\mathcal{A}}}}_{s}
P
\right)
\right]_{s}
\tag{1}
\\&=
\left[
\mathrm{CH}^{(1)}_{s}\left(
\mathrm{ip}^{(1,X)@}_{s}\left(
\mathrm{sc}^{0\mathbf{PT}_{\boldsymbol{\mathcal{A}}}}_{s}\left(
Q
\circ^{0\mathbf{PT}_{\boldsymbol{\mathcal{A}}}}_{s}
P
\right)
\right)
\right)
\right]_{s}
\tag{2}
\\&=
\left[
\mathrm{CH}^{(1)}_{s}\left(
\mathrm{sc}^{0\mathbf{Pth}_{\boldsymbol{\mathcal{A}}}}_{s}\left(
\left(
\mathrm{ip}^{(1,X)@}_{s}\left(
Q
\right)\right)
\circ^{0\mathbf{Pth}_{\boldsymbol{\mathcal{A}}}}_{s}
\left(
\mathrm{ip}^{(1,X)@}_{s}\left(
P
\right)\right)
\right)
\right)
\right]_{s}
\tag{3}
\\&=
\left[
\mathrm{CH}^{(1)}_{s}\left(
\mathrm{sc}^{0\mathbf{Pth}_{\boldsymbol{\mathcal{A}}}}_{s}\left(
\mathrm{ip}^{(1,X)@}_{s}\left(
P
\right)
\right)
\right)
\right]_{s}
\tag{4}
\\&=
\left[
\mathrm{CH}^{(1)}_{s}\left(
\mathrm{ip}^{(1,X)@}_{s}\left(
\mathrm{sc}^{0\mathbf{PT}_{\boldsymbol{\mathcal{A}}}}_{s}\left(
P
\right)
\right)
\right)
\right]_{s}
\tag{5}
\\&=
\left[
\mathrm{sc}^{0\mathbf{PT}_{\boldsymbol{\mathcal{A}}}}_{s}\left(
P
\right)
\right]_{s}
\tag{6}
\\&=
\mathrm{sc}^{0
[\mathbf{PT}_{\boldsymbol{\mathcal{A}}}]
}_{s}\left(
[P]_{s}
\right).
\tag{7}
\end{align*}
\end{flushleft}

The first equality unravels the interpretation of the $0$-source and $0$-composition operations in the many-sorted partial $\Sigma^{\boldsymbol{\mathcal{A}}}$-algebra $[\mathbf{PT}_{\boldsymbol{\mathcal{A}}}]$ according to Proposition~\ref{PPTQCatAlg}; the second equality follows from Lemma~\ref{LWCong}; the third equality follows from the fact that $\mathrm{ip}^{(1,X)@}$ is a $\Sigma^{\boldsymbol{\mathcal{A}}}$-homomorphism according to Definition~\ref{DIp}; the fourth equality follows from Proposition~\ref{PCHVarA4}; the fifth equality follows from the fact that $\mathrm{ip}^{(1,X)@}$ is a $\Sigma^{\boldsymbol{\mathcal{A}}}$-homomorphism according to Definition~\ref{DIp}; the sixth equality follows from Lemma~\ref{LWCong}; finally, the last equality recovers the interpretation of the  $0$-source and $0$-composition operations in the many-sorted partial $\Sigma^{\boldsymbol{\mathcal{A}}}$-algebra $[\mathbf{PT}_{\boldsymbol{\mathcal{A}}}]$ according to Proposition~\ref{PPTQCatAlg}.

The other equality holds by a similar argument.

This completes the proof.
\end{proof}

\begin{proposition}\label{PPTQVarA5} Let $s$ be a sort in $S$ and $[P]_{s}$ a path term class in $[\mathrm{PT}_{\boldsymbol{\mathcal{A}}}]_{s}$. Then the following equalities hold
\allowdisplaybreaks
\begin{align*}
[P]_{s}
\circ^{0
[\mathbf{PT}_{\boldsymbol{\mathcal{A}}}]
}_{s}
\left(
\mathrm{sc}^{0
[\mathbf{PT}_{\boldsymbol{\mathcal{A}}}]
}_{s}\left(
[P]_{s}
\right)
\right)
&=
[P]_{s};
\\
\left(
\mathrm{tg}^{0
[\mathbf{PT}_{\boldsymbol{\mathcal{A}}}]
}_{s}\left(
[P]_{s}
\right)
\right)
\circ^{0
[\mathbf{PT}_{\boldsymbol{\mathcal{A}}}]
}_{s}
[P]_{s}
&=
[P]_{s}.
\end{align*}
\end{proposition}
\begin{proof}
The following chain of equalities holds
\begin{flushleft}
$[P]_{s}
\circ^{0
[\mathbf{PT}_{\boldsymbol{\mathcal{A}}}]
}_{s}
\left(
\mathrm{sc}^{0
[\mathbf{PT}_{\boldsymbol{\mathcal{A}}}]
}_{s}\left(
[P]_{s}
\right)
\right)$
\allowdisplaybreaks
\begin{align*}
\qquad&=
\left[
P
\circ^{0\mathbf{PT}_{\boldsymbol{\mathcal{A}}}}_{s}
\left(
\mathrm{sc}^{0\mathbf{PT}_{\boldsymbol{\mathcal{A}}}}_{s}\left(
P
\right)
\right)
\right]_{s}
\tag{1}
\\&=
\left[
\mathrm{CH}^{(1)}_{s}\left(
\mathrm{ip}^{(1,X)@}_{s}\left(
P
\circ^{0\mathbf{PT}_{\boldsymbol{\mathcal{A}}}}_{s}
\left(
\mathrm{sc}^{0\mathbf{PT}_{\boldsymbol{\mathcal{A}}}}_{s}\left(
P
\right)
\right)
\right)
\right)
\right]_{s}
\tag{2}
\\&=
\left[
\mathrm{CH}^{(1)}_{s}\left(
\mathrm{ip}^{(1,X)@}_{s}\left(
P
\right)
\circ^{0\mathbf{Pth}_{\boldsymbol{\mathcal{A}}}}_{s}
\left(
\mathrm{sc}^{0\mathbf{Pth}_{\boldsymbol{\mathcal{A}}}}_{s}\left(
\mathrm{ip}^{(1,X)@}_{s}\left(
P
\right)
\right)
\right)
\right)
\right]_{s}
\tag{3}
\\&=
\left[
\mathrm{CH}^{(1)}_{s}\left(
\mathrm{ip}^{(1,X)@}_{s}\left(
P
\right)
\right)
\right]_{s}
\tag{4}
\\&=
[P]_{s}.
\tag{5}
\end{align*}
\end{flushleft}

The first equality unravels the interpretation of the $0$-source and $0$-composition operations in the many-sorted partial $\Sigma^{\boldsymbol{\mathcal{A}}}$-algebra $[\mathbf{PT}_{\boldsymbol{\mathcal{A}}}]$ according to Proposition~\ref{PPTQCatAlg}; the second equality follows from Lemma~\ref{LWCong}; the third equality follows from the fact that $\mathrm{ip}^{(1,X)@}$ is a $\Sigma^{\boldsymbol{\mathcal{A}}}$-homomorphism according to Definition~\ref{DIp}; the fourth equality follows from Proposition~\ref{PCHVarA5}; finally, the last equality follows from Lemma~\ref{LWCong}.

The other equality holds by a similar argument.

This completes the proof.
\end{proof}

\begin{proposition}\label{PPTQVarA6} Let $s$ be a sort in $S$ and $[P]_{s}$, $[Q]_{s}$ and $[R]_{s}$ path term classes in $[\mathrm{PT}_{\boldsymbol{\mathcal{A}}}]_{s}$ such that
\allowdisplaybreaks
\begin{align*}
\mathrm{sc}^{0
[\mathbf{PT}_{\boldsymbol{\mathcal{A}}}]
}_{s}\left(
[R]_{s}
\right)
&=
\mathrm{tg}^{0
[\mathbf{PT}_{\boldsymbol{\mathcal{A}}}]
}_{s}\left(
[Q]_{s}
\right);
\\
\mathrm{sc}^{0
[\mathbf{PT}_{\boldsymbol{\mathcal{A}}}]
}_{s}\left(
[Q]_{s}
\right)
&=
\mathrm{tg}^{0
[\mathbf{PT}_{\boldsymbol{\mathcal{A}}}]
}_{s}\left(
[P]_{s}
\right).
\end{align*}
Then the following equality holds
\[
[R]_{s}
\circ^{0
[\mathbf{PT}_{\boldsymbol{\mathcal{A}}}]
}_{s}
\left(
[Q]_{s}
\circ^{0
[\mathbf{PT}_{\boldsymbol{\mathcal{A}}}]
}_{s}
[P]_{s}
\right)
=
\left(
[R]_{s}
\circ^{0
[\mathbf{PT}_{\boldsymbol{\mathcal{A}}}]
}_{s}
[Q]_{s}
\right)
\circ^{0
[\mathbf{PT}_{\boldsymbol{\mathcal{A}}}]
}_{s}
[P]_{s}.
\]
\end{proposition}
\begin{proof}
The following chain of equalities holds
\begin{flushleft}
$[R]_{s}
\circ^{0
[\mathbf{PT}_{\boldsymbol{\mathcal{A}}}]
}_{s}
\left(
[Q]_{s}
\circ^{0
[\mathbf{PT}_{\boldsymbol{\mathcal{A}}}]
}_{s}
[P]_{s}
\right)$
\allowdisplaybreaks
\begin{align*}
&=
\left[
R
\circ^{0\mathbf{PT}_{\boldsymbol{\mathcal{A}}}}_{s}
\left(
Q
\circ^{0\mathbf{PT}_{\boldsymbol{\mathcal{A}}}}_{s}
P
\right)
\right]_{s}
\tag{1}
\\&=
\left[
\mathrm{CH}^{(1)}_{s}\left(
\mathrm{ip}^{(1,X)@}_{s}\left(
R
\circ^{0\mathbf{PT}_{\boldsymbol{\mathcal{A}}}}_{s}
\left(
Q
\circ^{0\mathbf{PT}_{\boldsymbol{\mathcal{A}}}}_{s}
P
\right)
\right)
\right)
\right]_{s}
\tag{2}
\\&=
\left[
\mathrm{CH}^{(1)}_{s}\left(
\mathrm{ip}^{(1,X)@}_{s}\left(
R
\right)
\circ^{0\mathbf{Pth}_{\boldsymbol{\mathcal{A}}}}_{s}
\left(
\mathrm{ip}^{(1,X)@}_{s}\left(
Q
\right)
\circ^{0\mathbf{Pth}_{\boldsymbol{\mathcal{A}}}}_{s}
\mathrm{ip}^{(1,X)@}_{s}\left(
P
\right)
\right)
\right)
\right]_{s}
\tag{3}
\\&=
\left[
\mathrm{CH}^{(1)}_{s}\left(
\left(
\mathrm{ip}^{(1,X)@}_{s}\left(
R
\right)
\circ^{0\mathbf{Pth}_{\boldsymbol{\mathcal{A}}}}_{s}
\mathrm{ip}^{(1,X)@}_{s}\left(
Q
\right)
\right)
\circ^{0\mathbf{Pth}_{\boldsymbol{\mathcal{A}}}}_{s}
\mathrm{ip}^{(1,X)@}_{s}\left(
P
\right)
\right)
\right]_{s}
\tag{4}
\\&=
\left[
\mathrm{CH}^{(1)}_{s}\left(
\mathrm{ip}^{(1,X)@}_{s}\left(
\left(
R
\circ^{0
[\mathbf{PT}_{\boldsymbol{\mathcal{A}}}]}_{s}
Q
\right)
\circ^{0
[\mathbf{PT}_{\boldsymbol{\mathcal{A}}}]}_{s}
P
\right)\right)
\right]_{s}
\tag{5}
\\&=
\left[
\left(
R
\circ^{0
[\mathbf{PT}_{\boldsymbol{\mathcal{A}}}]}_{s}
Q
\right)
\circ^{0
[\mathbf{PT}_{\boldsymbol{\mathcal{A}}}]}_{s}
P
\right]_{s}
\tag{6}
\\&=
\left(
[R]_{s}
\circ^{0
[\mathbf{PT}_{\boldsymbol{\mathcal{A}}}]
}_{s}
[Q]_{s}
\right)
\circ^{0
[\mathbf{PT}_{\boldsymbol{\mathcal{A}}}]
}_{s}
[P]_{s}.
\tag{7}
\end{align*}
\end{flushleft}

The first equality unravels the interpretation of the $0$-composition operation in the many-sorted partial $\Sigma^{\boldsymbol{\mathcal{A}}}$-algebra $[\mathbf{PT}_{\boldsymbol{\mathcal{A}}}]$ according to Proposition~\ref{PPTQCatAlg}; the second equality follows from Lemma~\ref{LWCong}; the third equality follows from the fact that $\mathrm{ip}^{(1,X)@}$ is a $\Sigma^{\boldsymbol{\mathcal{A}}}$-homomorphism according to Definition~\ref{DIp}; the fourth equality follows from Proposition~\ref{PCHVarA6}; the fifth equality follows from the fact that $\mathrm{ip}^{(1,X)@}$ is a $\Sigma^{\boldsymbol{\mathcal{A}}}$-homomorphism according to Definition~\ref{DIp}; the sixth equality follows from Lemma~\ref{LWCong}; finally, the last equality recovers the interpretation of the  $0$-composition operation in the many-sorted partial $\Sigma^{\boldsymbol{\mathcal{A}}}$-algebra $[\mathbf{PT}_{\boldsymbol{\mathcal{A}}}]$ according to Proposition~\ref{PPTQCatAlg}.

This completes the proof.
\end{proof}

On the basis of the above propositions, we can now show that $[\mathrm{PT}_{\boldsymbol{\mathcal{A}}}]$ is equipped with a structure of $S$-sorted category.

\begin{restatable}{definition}{DPTQCat}
\label{DPTQCat}
\index{path terms!first-order!$[\mathsf{PT}_{\boldsymbol{\mathcal{A}}}]$}
Let $[\mathsf{PT}_{\boldsymbol{\mathcal{A}}}]$ stand for the ordered pair
\[
[\mathsf{PT}_{\boldsymbol{\mathcal{A}}}]
=
\left(
[\mathrm{PT}_{\boldsymbol{\mathcal{A}}}],
\left(
\circ^{0[\mathbf{PT}_{\boldsymbol{\mathcal{A}}}]},
\mathrm{sc}^{0[\mathbf{PT}_{\boldsymbol{\mathcal{A}}}]},
\mathrm{tg}^{0[\mathbf{PT}_{\boldsymbol{\mathcal{A}}}]}
\right)
\right).
\]
\end{restatable}

\begin{restatable}{proposition}{PPTQCat}
\label{PPTQCat} $[\mathsf{PT}_{\boldsymbol{\mathcal{A}}}]$ is an $S$-sorted category.
\end{restatable}
\begin{proof}
Let us recall that, for every $s\in S$, the structure
$([\mathrm{PT}_{\boldsymbol{\mathcal{A}}}]_{s},\xi_{0,s})$, where
$$
\xi_{0,s}=
\left(
\circ^{0
[\mathbf{PT}_{\boldsymbol{\mathcal{A}}}]
}_{s},
\mathrm{sc}^{0
[\mathbf{PT}_{\boldsymbol{\mathcal{A}}}]
}_{s},
\mathrm{tg}^{0
[\mathbf{PT}_{\boldsymbol{\mathcal{A}}}]
}_{s}
\right)
$$
is a single-sorted category in virtue of Definition~\ref{DCat} and Propositions~\ref{PPTQVarA2}, \ref{PPTQVarA3}, \ref{PPTQVarA4}, \ref{PPTQVarA5} and~\ref{PPTQVarA6}.
Therefore, following Definition~\ref{DnCat}, the $S$-sorted structure
$$\left(
[\mathrm{PT}_{\boldsymbol{\mathcal{A}}}]_{s},\xi_{0,s}
\right)_{s\in S}$$
is an $S$-sorted category.
\end{proof}


Our next aim is to prove that $[\mathsf{PT}_{\boldsymbol{\mathcal{A}}}]$ is an $S$-sorted categorial $\Sigma$-algebra. To state this fact properly we prove the following propositions.

\begin{proposition}\label{PPTQVarA7} For every $(\mathbf{s},s)\in S^{\star}\times S$, every operation symbol $\sigma\in \Sigma_{\mathbf{s},s}$ and every family of path term classes $([P_{j}]_{s_{j}})_{j\in\bb{\mathbf{s}}}$ in $[\mathbf{PT}_{\boldsymbol{\mathcal{A}}}]_{\mathbf{s}}$, the following equalities holds
\begin{align*}
\mathrm{sc}^{0
[\mathbf{PT}_{\boldsymbol{\mathcal{A}}}]}_{s}\left(
\sigma^{
[\mathbf{PT}_{\boldsymbol{\mathcal{A}}}]}
\left(\left(
\left[
P_{j}
\right]_{s_{j}}
\right)_{j\in\bb{\mathbf{s}}}
\right)\right)
&=
\sigma^{
[\mathbf{PT}_{\boldsymbol{\mathcal{A}}}]}
\left(\left(
\mathrm{sc}^{0
[\mathbf{PT}_{\boldsymbol{\mathcal{A}}}]}_{s_{j}}\left(
\left[
P_{j}
\right]_{s_{j}}
\right)
\right)_{j\in\bb{\mathbf{s}}}
\right);
\\
\mathrm{tg}^{0
[\mathbf{PT}_{\boldsymbol{\mathcal{A}}}]}_{s}\left(
\sigma^{
[\mathbf{PT}_{\boldsymbol{\mathcal{A}}}]}
\left(\left(
\left[
P_{j}
\right]_{s_{j}}
\right)_{j\in\bb{\mathbf{s}}}
\right)\right)
&=
\sigma^{
[\mathbf{PT}_{\boldsymbol{\mathcal{A}}}]}
\left(\left(
\mathrm{tg}^{0
[\mathbf{PT}_{\boldsymbol{\mathcal{A}}}]}_{s_{j}}\left(
\left[
P_{j}
\right]_{s_{j}}
\right)
\right)_{j\in\bb{\mathbf{s}}}
\right).
\end{align*}
\end{proposition}
\begin{proof}
The following chain of equalities holds
\begin{flushleft}
$\mathrm{sc}^{0
[\mathbf{PT}_{\boldsymbol{\mathcal{A}}}]}_{s}\left(
\sigma^{
[\mathbf{PT}_{\boldsymbol{\mathcal{A}}}]}
\left(\left(
\left[
P_{j}
\right]_{s_{j}}
\right)_{j\in\bb{\mathbf{s}}}
\right)\right)$
\allowdisplaybreaks
\begin{align*}
\qquad&=
\left[
\mathrm{sc}^{0\mathbf{PT}_{\boldsymbol{\mathcal{A}}}}_{s}\left(
\sigma^{\mathbf{PT}_{\boldsymbol{\mathcal{A}}}}
\left(\left(
P_{j}
\right)_{j\in\bb{\mathbf{s}}}
\right)
\right)
\right]_{s}
\tag{1}
\\&=
\left[
\mathrm{CH}^{(1)}_{s}\left(
\mathrm{ip}^{(1,X)@}_{s}\left(
\mathrm{sc}^{0\mathbf{PT}_{\boldsymbol{\mathcal{A}}}}_{s}\left(
\sigma^{\mathbf{PT}_{\boldsymbol{\mathcal{A}}}}
\left(\left(
P_{j}
\right)_{j\in\bb{\mathbf{s}}}
\right)
\right)
\right)
\right)
\right]_{s}
\tag{2}
\\&=
\left[
\mathrm{CH}^{(1)}_{s}\left(
\mathrm{sc}^{0\mathbf{Pth}_{\boldsymbol{\mathcal{A}}}}_{s}\left(
\sigma^{\mathbf{Pth}_{\boldsymbol{\mathcal{A}}}}
\left(\left(
\mathrm{ip}^{(1,X)@}_{s_{j}}\left(
P_{j}
\right)
\right)_{j\in\bb{\mathbf{s}}}
\right)
\right)
\right)
\right]_{s}
\tag{3}
\\&=
\left[
\mathrm{CH}^{(1)}_{s}\left(
\sigma^{\mathbf{Pth}_{\boldsymbol{\mathcal{A}}}}
\left(\left(
\mathrm{sc}^{0\mathbf{Pth}_{\boldsymbol{\mathcal{A}}}}_{s_{j}}\left(
\mathrm{ip}^{(1,X)@}_{s_{j}}\left(
P_{j}
\right)
\right)
\right)_{j\in\bb{\mathbf{s}}}
\right)
\right)
\right]_{s}
\tag{4}
\\&=
\left[
\mathrm{CH}^{(1)}_{s}\left(
\mathrm{ip}^{(1,X)@}_{s}\left(
\sigma^{\mathbf{PT}_{\boldsymbol{\mathcal{A}}}}
\left(\left(
\mathrm{sc}^{0\mathbf{PT}_{\boldsymbol{\mathcal{A}}}}_{s_{j}}\left(
P_{j}
\right)
\right)_{j\in\bb{\mathbf{s}}}
\right)
\right)
\right)
\right]_{s}
\tag{5}
\\&=
\left[
\sigma^{\mathbf{PT}_{\boldsymbol{\mathcal{A}}}}
\left(\left(
\mathrm{sc}^{0\mathbf{PT}_{\boldsymbol{\mathcal{A}}}}_{s_{j}}\left(
P_{j}
\right)
\right)_{j\in\bb{\mathbf{s}}}
\right)
\right]_{s}
\tag{6}
\\&=
\sigma^{
[\mathbf{PT}_{\boldsymbol{\mathcal{A}}}]}
\left(\left(
\mathrm{sc}^{0
[\mathbf{PT}_{\boldsymbol{\mathcal{A}}}]}_{s_{j}}\left(
\left[
P_{j}
\right]_{s_{j}}
\right)
\right)_{j\in\bb{\mathbf{s}}}
\right).
\tag{7}
\end{align*}
\end{flushleft}

The first equality unravels the interpretation of the $\sigma$ and $0$-source operations in the many-sorted partial $\Sigma^{\boldsymbol{\mathcal{A}}}$-algebra $[\mathbf{PT}_{\boldsymbol{\mathcal{A}}}]$ according to Proposition~\ref{PPTQCatAlg}; the second equality follows from Lemma~\ref{LWCong}; the third equality follows from the fact that $\mathrm{ip}^{(1,X)@}$ is a $\Sigma^{\boldsymbol{\mathcal{A}}}$-homomorphism according to Definition~\ref{DIp}; the fourth equality follows from Proposition~\ref{PCHVarA7}; the fifth equality follows from the fact that $\mathrm{ip}^{(1,X)@}$ is a $\Sigma^{\boldsymbol{\mathcal{A}}}$-homomorphism according to Definition~\ref{DIp}; the sixth equality follows from Lemma~\ref{LWCong}; finally, the last equality recovers the interpretation of the $\sigma$ and $0$-source operations in the many-sorted partial $\Sigma^{\boldsymbol{\mathcal{A}}}$-algebra $[\mathbf{PT}_{\boldsymbol{\mathcal{A}}}]$ according to Proposition~\ref{PPTQCatAlg}.

The other equality holds by a similar argument.

This completes the proof.
\end{proof}

\begin{proposition}\label{PPTQVarA8} Let $(\mathbf{s},s)$ be an element in $S^{\star}\times S$, $\sigma$ an operation symbol in $\Sigma_{\mathbf{s},s}$ and $([P_{j}]_{s_{j}})_{j\in\bb{\mathbf{s}}}$, $([Q_{j}]_{s_{j}})_{j\in\bb{\mathbf{s}}}$ families of path term classes in $[\mathrm{PT}_{\boldsymbol{\mathcal{A}}}]_{\mathbf{s}}$ such that, for every $j\in\bb{\mathbf{s}}$, 
$$
\mathrm{sc}^{0
[\mathbf{PT}_{\boldsymbol{\mathcal{A}}}]
}_{s_{j}}\left(
[Q_{j}]_{s_{j}}
\right)
=
\mathrm{tg}^{0
[\mathbf{PT}_{\boldsymbol{\mathcal{A}}}]
}_{s_{j}}\left(
[P_{j}]_{s_{j}}
\right).
$$
Then the following chain of equalities holds
\begin{multline*}
\sigma^{
[\mathbf{PT}_{\boldsymbol{\mathcal{A}}}]}
\left(\left(
[Q_{j}]_{s_{j}}
\circ^{0
[\mathbf{PT}_{\boldsymbol{\mathcal{A}}}]}_{s_{j}}
[P_{j}]_{s_{j}}
\right)_{j\in\bb{\mathbf{s}}}
\right)
\\
=
\sigma^{
[\mathbf{PT}_{\boldsymbol{\mathcal{A}}}]}
\left(\left(
[Q_{j}]_{s_{j}}
\right)_{j\in\bb{\mathbf{s}}}
\right)
\circ^{0
[\mathbf{PT}_{\boldsymbol{\mathcal{A}}}]}_{s}
\sigma^{
[\mathbf{PT}_{\boldsymbol{\mathcal{A}}}]}
\left(\left(
[P_{j}]_{s_{j}}
\right)_{j\in\bb{\mathbf{s}}}
\right).
\end{multline*}
\end{proposition}
\begin{proof}
The following chain of equalities holds
\begin{flushleft}
$\sigma^{
[\mathbf{PT}_{\boldsymbol{\mathcal{A}}}]}
\left(\left(
[Q_{j}]_{s_{j}}
\circ^{0
[\mathbf{PT}_{\boldsymbol{\mathcal{A}}}]}_{s_{j}}
[P_{j}]_{s_{j}}
\right)_{j\in\bb{\mathbf{s}}}
\right)$
\allowdisplaybreaks
\begin{align*}
&=
\left[
\sigma^{\mathbf{PT}_{\boldsymbol{\mathcal{A}}}}
\left(\left(
Q_{j}
\circ^{0\mathbf{PT}_{\boldsymbol{\mathcal{A}}}}_{s_{j}}
P_{j}
\right)_{j\in\bb{\mathbf{s}}}
\right)
\right]_{s}
\tag{1}
\\&=
\left[
\mathrm{CH}^{(1)}_{s}\left(
\mathrm{ip}^{(1,X)@}_{s}\left(
\sigma^{\mathbf{PT}_{\boldsymbol{\mathcal{A}}}}
\left(\left(
Q_{j}
\circ^{0\mathbf{PT}_{\boldsymbol{\mathcal{A}}}}_{s_{j}}
P_{j}
\right)_{j\in\bb{\mathbf{s}}}
\right)
\right)\right)
\right]_{s}
\tag{2}
\\&=
\left[
\mathrm{CH}^{(1)}_{s}\left(
\sigma^{\mathbf{Pth}_{\boldsymbol{\mathcal{A}}}}
\left(\left(
\mathrm{ip}^{(1,X)@}_{s_{j}}\left(
Q_{j}
\right)
\circ^{0\mathbf{Pth}_{\boldsymbol{\mathcal{A}}}}_{s_{j}}
\mathrm{ip}^{(1,X)@}_{s_{j}}\left(
P_{j}
\right)
\right)_{j\in\bb{\mathbf{s}}}\right)
\right)
\right]_{s}
\tag{3}
\\&=
\left[
\mathrm{CH}^{(1)}_{s}\left(
\sigma^{\mathbf{Pth}_{\boldsymbol{\mathcal{A}}}}
\left(\left(
\mathrm{ip}^{(1,X)@}_{s_{j}}\left(
Q_{j}
\right)\right)_{j\in\bb{\mathbf{s}}}\right)
\circ^{0\mathbf{Pth}_{\boldsymbol{\mathcal{A}}}}_{s}
\right.\right.
\\&\qquad\qquad\qquad\qquad\qquad\qquad\qquad\qquad\qquad
\left.\left.
\sigma^{\mathbf{Pth}_{\boldsymbol{\mathcal{A}}}}
\left(\left(
\mathrm{ip}^{(1,X)@}_{s_{j}}\left(
P_{j}
\right)\right)_{j\in\bb{\mathbf{s}}}\right)
\right)
\right]_{s}
\tag{4}
\\&=
\left[
\mathrm{CH}^{(1)}_{s}\left(
\mathrm{ip}^{(1,X)@}_{s}\left(
\sigma^{\mathbf{PT}_{\boldsymbol{\mathcal{A}}}}
\left(\left(
Q_{j}
\right)_{j\in\bb{\mathbf{s}}}\right)
\circ^{0\mathbf{PT}_{\boldsymbol{\mathcal{A}}}}_{s}
\sigma^{\mathbf{PT}_{\boldsymbol{\mathcal{A}}}}
\left(\left(
P_{j}
\right)_{j\in\bb{\mathbf{s}}}\right)
\right)\right)
\right]_{s}
\tag{5}
\\&=
\left[
\sigma^{\mathbf{PT}_{\boldsymbol{\mathcal{A}}}}
\left(\left(
Q_{j}
\right)_{j\in\bb{\mathbf{s}}}\right)
\circ^{0\mathbf{PT}_{\boldsymbol{\mathcal{A}}}}_{s}
\sigma^{\mathbf{PT}_{\boldsymbol{\mathcal{A}}}}
\left(\left(
P_{j}
\right)_{j\in\bb{\mathbf{s}}}\right)
\right]_{s}
\tag{6}
\\&=
\sigma^{
[\mathbf{PT}_{\boldsymbol{\mathcal{A}}}]}
\left(\left(
[Q_{j}]_{s_{j}}
\right)_{j\in\bb{\mathbf{s}}}
\right)
\circ^{0
[\mathbf{PT}_{\boldsymbol{\mathcal{A}}}]}_{s}
\sigma^{
[\mathbf{PT}_{\boldsymbol{\mathcal{A}}}]}
\left(\left(
[P_{j}]_{s_{j}}
\right)_{j\in\bb{\mathbf{s}}}
\right).
\tag{7}
\end{align*}
\end{flushleft}

The first equality unpacks the interpretation of the $\sigma$ and $0$-composition operations in the many-sorted partial $\Sigma^{\boldsymbol{\mathcal{A}}}$-algebra $[\mathbf{PT}_{\boldsymbol{\mathcal{A}}}]$ according to Proposition~\ref{PPTQCatAlg}; the second equality follows from Lemma~\ref{LWCong}; the third equality follows from the fact that $\mathrm{ip}^{(1,X)@}$ is a $\Sigma^{\boldsymbol{\mathcal{A}}}$-homomorphism according to Definition~\ref{DIp}; the fourth equality follows from Proposition~\ref{PCHVarA8}; the fifth equality follows from the fact that $\mathrm{ip}^{(1,X)@}$ is a $\Sigma^{\boldsymbol{\mathcal{A}}}$-homomorphism according to Definition~\ref{DIp}; the sixth equality follows from Lemma~\ref{LWCong}; finally, the last equality recovers the interpretation of the $\sigma$ and   $0$-composition operations in the many-sorted partial $\Sigma^{\boldsymbol{\mathcal{A}}}$-algebra $[\mathbf{PT}_{\boldsymbol{\mathcal{A}}}]$ according to Proposition~\ref{PPTQCatAlg}.

This completes the proof.
\end{proof}

\begin{restatable}{proposition}{PPTQCtyAlg}\label{PPTQCtyAlg}
$[\mathsf{PT}_{\boldsymbol{\mathcal{A}}}]$
is an $S$-sorted categorial $\Sigma$-algebra.
\end{restatable}
\begin{proof}
That $[\mathsf{PT}_{\boldsymbol{\mathcal{A}}}]$ is an $S$-sorted category was already proven in Proposition~\ref{PPTQCat}. We are only left to prove that, for every $(\mathbf{s},s)$ in $S^{\star}\times S$, $\sigma^{[\mathsf{PT}_{\boldsymbol{\mathcal{A}}}(X)]}$ is a functor from the product category $[\mathsf{PT}_{\boldsymbol{\mathcal{A}}}]_{\mathbf{s}}$ to the category $[\mathsf{PT}_{\boldsymbol{\mathcal{A}}}]_{s}$.  This is a fact that follows from Propositions~\ref{PPTQVarA7} and~\ref{PPTQVarA8}.

This completes the proof.
\end{proof}

\section{
\texorpdfstring
{An Artinian order on $\coprod[\mathrm{PT}_{\boldsymbol{\mathcal{A}}}]$}
{An Artinian order on the quotient of path terms}
}
In this section we define on $\coprod[\mathrm{PT}_{\boldsymbol{\mathcal{A}}}]$ an Artinian  order. 

\begin{restatable}{definition}{DPTQOrd}
\label{DPTQOrd} 
\index{partial order!first-order!$\leq_{[\mathbf{PT}_{\boldsymbol{\mathcal{A}}}]}$}
We let $\leq_{[\mathbf{PT}_{\boldsymbol{\mathcal{A}}}]}$ stand for the binary relation on 
$\coprod [\mathrm{PT}_{\boldsymbol{\mathcal{A}}}]$ which consists of those ordered pairs
$(([Q]_{t},t),([P]_{s},s))$ in $(\coprod [\mathrm{PT}_{\boldsymbol{\mathcal{A}}}])^{2}$ for which there exists $Q'\in [Q]_{t}$ and $P'\in [P]_{s}$ such that
$$
\left(
\mathrm{ip}^{(1,X)@}_{t}\left(
Q'
\right),t\right)
\leq_{\mathbf{Pth}_{\boldsymbol{\mathcal{A}}}}
\left(
\mathrm{ip}^{(1,X)@}_{s}\left(
P'
\right),s
\right)
.
$$
\end{restatable}

We now prove that $\leq_{[\mathbf{PT}_{\boldsymbol{\mathcal{A}}}]}$  is an Artinian order on 
$\coprod [\mathrm{PT}_{\boldsymbol{\mathcal{A}}}]$.

\begin{restatable}{proposition}{PPTQOrdArt}
\label{PPTQOrdArt} $(\coprod[\mathrm{PT}_{\boldsymbol{\mathcal{A}}}], \leq_{[\mathbf{PT}_{\boldsymbol{\mathcal{A}}}]})$ is an Artinian ordered set. 
\end{restatable}
\begin{proof}
We first prove that $\leq_{[\mathbf{PT}_{\boldsymbol{\mathcal{A}}}]}$ is reflexive. Note that for every sort $s\in S$ and every path term $P\in\mathrm{PT}_{\boldsymbol{\mathcal{A}},s}$, the inequality
\[
\left(\mathrm{ip}^{(1,X)@}_{s}(P),s\right)
\leq_{\mathbf{Pth}_{\boldsymbol{\mathcal{A}}}}
\left(\mathrm{ip}^{(1,X)@}_{s}(P),s\right)
\]
holds since $\leq_{\mathbf{Pth}_{\boldsymbol{\mathcal{A}}}}$ is a reflexive relation on $\coprod\mathrm{Pth}_{\boldsymbol{\mathcal{A}}}$, as seen in Definition~\ref{DOrd}. Thus, $\leq_{[\mathbf{PT}_{\boldsymbol{\mathcal{A}}}]}$ is reflexive.

We next prove that $\leq_{[\mathbf{PT}_{\boldsymbol{\mathcal{A}}}]}$ is antisymmetric. Let $s, t$ be sorts in $S$ and $[P]_{s}, [Q]_{t}$ be path term classes in $[\mathrm{PT}_{\boldsymbol{\mathcal{A}}}]_{s}$ and $[\mathrm{PT}_{\boldsymbol{\mathcal{A}}}]_{t}$, respectively, satisfying that 
\begin{align*}
\left([P]_{s},s\right)
&\leq_{[\mathbf{PT}_{\boldsymbol{\mathcal{A}}}]}
\left([Q]_{t},t\right);
&
\left([Q]_{t},t\right)
&\leq_{[\mathbf{PT}_{\boldsymbol{\mathcal{A}}}]}
\left([P]_{s},s\right).
\end{align*}

Following Definition~\ref{DPTQOrd}, there exists path term representatives, $P',P''\in [P]_{s}$ and $Q',Q''\in [Q]_{t}$ satisfying that 
\begin{align*}
\left(\mathrm{ip}^{(1,X)@}_{s}(P'),s\right)
&\leq_{\mathbf{Pth}_{\boldsymbol{\mathcal{A}}}}
\left(\mathrm{ip}^{(1,X)@}_{t}(Q'),t\right);
\\
\left(\mathrm{ip}^{(1,X)@}_{t}(Q''),t\right)
&\leq_{\mathbf{Pth}_{\boldsymbol{\mathcal{A}}}}
\left(\mathrm{ip}^{(1,X)@}_{s}(P''),s\right).
\end{align*}

Following Proposition~\ref{PCHMono}, we conclude that
\begin{align*}
\left(\mathrm{CH}^{(1)}_{s}\left(\mathrm{ip}^{(1,X)@}_{s}(P')\right),s\right)
&\leq_{\mathbf{T}_{\Sigma^{\boldsymbol{\mathcal{A}}}}(X)}
\left(\mathrm{CH}^{(1)}_{t}\left(\mathrm{ip}^{(1,X)@}_{t}(Q')\right),t\right);
\\
\left(\mathrm{CH}^{(1)}_{t}\left(\mathrm{ip}^{(1,X)@}_{t}(Q'')\right),t\right)
&\leq_{\mathbf{T}_{\Sigma^{\boldsymbol{\mathcal{A}}}}(X)}
\left(\mathrm{CH}^{(1)}_{s}\left(\mathrm{ip}^{(1,X)@}_{s}(P'')\right),s\right);
\end{align*}

Note that by Lemma~\ref{LThetaCong}, we have that 
\begin{align*}
\mathrm{CH}^{(1)}_{s}\left(
\mathrm{ip}^{(1,X)@}_{s}\left(
P'
\right)\right)&=
\mathrm{CH}^{(1)}_{s}\left(
\mathrm{ip}^{(1,X)@}_{s}\left(
P''
\right)\right);
\\
\mathrm{CH}^{(1)}_{t}\left(
\mathrm{ip}^{(1,X)@}_{t}\left(
Q'
\right)\right)&=
\mathrm{CH}^{(1)}_{t}\left(
\mathrm{ip}^{(1,X)@}_{t}\left(
Q''
\right)\right);
\end{align*}

Therefore, since $\leq_{\mathbf{T}_{\Sigma^{\boldsymbol{\mathcal{A}}}}(X)}$ is a partial order on $\coprod\mathrm{T}_{\Sigma^{\boldsymbol{\mathcal{A}}}}(X)$, we conclude that $s=t$ and 
\[
\mathrm{CH}^{(1)}_{s}\left(
\mathrm{ip}^{(1,X)@}_{s}\left(
P
\right)\right)=
\mathrm{CH}^{(1)}_{s}\left(
\mathrm{ip}^{(1,X)@}_{s}\left(
Q
\right)\right).
\]

Now, following Lemma~\ref{LWCong}, we have that 
\begin{align*}
\left(
P,
\mathrm{CH}^{(1)}_{s}\left(
\mathrm{ip}^{(1,X)@}_{s}\left(
P
\right)\right)
\right)
&\in\Theta^{[1]}_{s};
&
\left(
Q,
\mathrm{CH}^{(1)}_{s}\left(
\mathrm{ip}^{(1,X)@}_{s}\left(
Q
\right)\right)
\right)
&\in\Theta^{[1]}_{s}.
\end{align*}

Since $\Theta^{[1]}$ is an equivalence relation, we conclude that $[P]_{s}=[Q]_{s}$. Thus, $\leq_{[\mathbf{PT}_{\boldsymbol{\mathcal{A}}}]}$ is antisymmetric.

We next prove that $\leq_{[\mathbf{PT}_{\boldsymbol{\mathcal{A}}}]}$ is transitive. Let $s, t, u$ be sorts in $S$ and $[P]_{s}, [Q]_{t}, [R]_{u}$ be path term classes in $[\mathrm{PT}_{\boldsymbol{\mathcal{A}},s}]_{s}$,$ [\mathrm{PT}_{\boldsymbol{\mathcal{A}},t}]_{t}$  and $[\mathrm{PT}_{\boldsymbol{\mathcal{A}},u}]_{u}$, respectively, satisfying that 
\begin{align*}
\left([P]_{s},s\right)
&\leq_{[\mathbf{PT}_{\boldsymbol{\mathcal{A}}}]}
\left([Q]_{t},t\right);
&
\left([Q]_{t},t\right)
&\leq_{[\mathbf{PT}_{\boldsymbol{\mathcal{A}}}]}
\left([R]_{u},u\right).
\end{align*}

Following Definition~\ref{DPTQOrd}, there exists path term representatives, $P'\in [P]_{s}$, $Q', Q''\in [Q]_{t}$ and $R'\in [R]_{u}$  such that 
\begin{align*}
\left(\mathrm{ip}^{(1,X)@}_{s}(P'),s\right)
&\leq_{\mathbf{Pth}_{\boldsymbol{\mathcal{A}}}}
\left(\mathrm{ip}^{(1,X)@}_{t}(Q'),t\right);
\\
\left(\mathrm{ip}^{(1,X)@}_{t}(Q''),t\right)
&\leq_{\mathbf{Pth}_{\boldsymbol{\mathcal{A}}}}
\left(\mathrm{ip}^{(1,X)@}_{u}(R),u\right).
\end{align*}

Following Proposition~\ref{PCHMono}, we conclude that
\begin{align*}
\left(\mathrm{CH}^{(1)}_{s}\left(\mathrm{ip}^{(1,X)@}_{s}(P')\right),s\right)
&\leq_{\mathbf{T}_{\Sigma^{\boldsymbol{\mathcal{A}}}}(X)}
\left(\mathrm{CH}^{(1)}_{t}\left(\mathrm{ip}^{(1,X)@}_{t}(Q')\right),t\right);
\\
\left(\mathrm{CH}^{(1)}_{t}\left(\mathrm{ip}^{(1,X)@}_{t}(Q'')\right),t\right)
&\leq_{\mathbf{T}_{\Sigma^{\boldsymbol{\mathcal{A}}}}(X)}
\left(\mathrm{CH}^{(1)}_{u}\left(\mathrm{ip}^{(1,X)@}_{u}(R')\right),u\right);
\end{align*}

Note that by Lemma~\ref{LThetaCong}, we have that 
\begin{align*}
\mathrm{CH}^{(1)}_{t}\left(
\mathrm{ip}^{(1,X)@}_{t}\left(
Q'
\right)\right)&=
\mathrm{CH}^{(1)}_{t}\left(
\mathrm{ip}^{(1,X)@}_{t}\left(
Q''
\right)\right).
\end{align*}

Therefore, since $\leq_{\mathbf{T}_{\Sigma^{\boldsymbol{\mathcal{A}}}}(X)}$ is an order on $\coprod\mathrm{T}_{\Sigma^{\boldsymbol{\mathcal{A}}}}(X)$, we conclude that
\begin{align*}
\left(\mathrm{CH}^{(1)}_{s}\left(\mathrm{ip}^{(1,X)@}_{s}(P')\right),s\right)
&\leq_{\mathbf{T}_{\Sigma^{\boldsymbol{\mathcal{A}}}}(X)}
\left(\mathrm{CH}^{(1)}_{u}\left(\mathrm{ip}^{(1,X)@}_{u}(R')\right),u\right).
\end{align*}

Now, by Proposition~\ref{PIpCHOrd}, the following inequality holds
\begin{multline*}
\left(\mathrm{ip}^{(1,X)@}_{s}\left(\mathrm{CH}^{(1)}_{s}\left(\mathrm{ip}^{(1,X)@}_{s}(P')\right)\right),s\right)
\\
\leq_{\mathbf{Pth}_{\boldsymbol{\mathcal{A}}}}
\left(\mathrm{ip}^{(1,X)@}_{u}\left(\mathrm{CH}^{(1)}_{u}\left(\mathrm{ip}^{(1,X)@}_{u}(R')\right)\right),u\right).
\tag{$\star$}
\end{multline*}

Note that by Lemma~\ref{LThetaCong}, we have that 
\allowdisplaybreaks
\begin{align*}
\mathrm{CH}^{(1)}_{s}\left(
\mathrm{ip}^{(1,X)@}_{s}\left(
P'
\right)\right)&=
\mathrm{CH}^{(1)}_{s}\left(
\mathrm{ip}^{(1,X)@}_{s}\left(
P
\right)\right);\\
\mathrm{CH}^{(1)}_{u}\left(
\mathrm{ip}^{(1,X)@}_{u}\left(
R'
\right)\right)&=
\mathrm{CH}^{(1)}_{u}\left(
\mathrm{ip}^{(1,X)@}_{u}\left(
R
\right)\right);
\end{align*}

Moreover, by Lemma~\ref{LWCong}, we have that	
\allowdisplaybreaks
\begin{align*}
\mathrm{CH}^{(1)}_{s}\left(
\mathrm{ip}^{(1,X)@}_{s}\left(
P
\right)\right)
&\in [P]_{s};
&
\mathrm{CH}^{(1)}_{u}\left(
\mathrm{ip}^{(1,X)@}_{u}\left(
R
\right)\right)
&\in [R]_{u}.
\end{align*}

Thus, taking into account the inequality $(\star)$ and Definition~\ref{DPTQOrd}, we conclude that the following inequality holds
\[
\left(
[P]_{s},s
\right)
\leq_{[\mathbf{PT}_{\boldsymbol{\mathcal{A}}}]}
\left(
[R]_{u},u
\right).
\]
Thus, $\leq_{[\mathbf{PT}_{\boldsymbol{\mathcal{A}}}]}$ is transitive.

It follows that $(\coprod[\mathrm{PT}_{\boldsymbol{\mathcal{A}}}], \leq_{[\mathbf{PT}_{\boldsymbol{\mathcal{A}}}]})$ is a partially ordered set.

We next prove that no strictly decreasing $\omega_{0}$-chain can exist in the ordered set 
$(\coprod[\mathrm{PT}_{\boldsymbol{\mathcal{A}}}], \leq_{[\mathbf{PT}_{\boldsymbol{\mathcal{A}}}]})$. 

Let us suppose, towards a contradiction, that there exists at least one strictly decreasing $\omega_{0}$-chain $(([R_{k}]_{s_{k}},s_{k}))_{k\in\omega_{0}}$ in $(\coprod[\mathrm{PT}_{\boldsymbol{\mathcal{A}}}], \leq_{[\mathbf{PT}_{\boldsymbol{\mathcal{A}}}]})$.

Then, for every $k\in\omega_{0}$, we have that 
\begin{enumerate}
\item[(i)] $[R_{k}]_{s_{k}}$ is a path term class in $[\mathrm{PT}_{\boldsymbol{\mathcal{A}}}]_{s_{k}}$, and 
\item[(ii)] $([R_{k+1}]_{s_{k+1}},s_{k+1})
<_{[\mathbf{PT}_{\boldsymbol{\mathcal{A}}}]}
([R_{k}]_{s_{k}},s_{k})$.
\end{enumerate}

By Definition~\ref{DPTQOrd}, for every $k\in\omega_{0}$, there exists path term representatives $R'$ and $R''$ in $[R_{k}]_{s_{k}}$ for which
\[
\left(
\mathrm{ip}^{(1,X)@}_{s_{k+1}}\left(R'_{k+1}\right), s_{k+1}
\right)
<_{\mathbf{Pth}_{\boldsymbol{\mathcal{A}}}}
\left(
\mathrm{ip}^{(1,X)@}_{s_{k}}\left(R''_{k}\right), s_{k}
\right).
\]

By Corollary~\ref{CIpCHOrd}, for every $k\in\omega_{0}$, the following inequality holds
\begin{multline*}
\left(\mathrm{ip}^{(1,X)@}_{s_{k+1}}\left(\mathrm{CH}^{(1)}_{s_{k+1}}\left(\mathrm{ip}^{(1,X)@}_{s_{k+1}}(R'_{k+1})\right)\right),s_{k+1}\right)
\\
<_{\mathbf{Pth}_{\boldsymbol{\mathcal{A}}}}
\left(\mathrm{ip}^{(1,X)@}_{s_{k}}\left(\mathrm{CH}^{(1)}_{s_{k}}\left(\mathrm{ip}^{(1,X)@}_{s_{k}}(R''_{s_{k}})\right)\right),s_{k}\right).
\tag{$\star$}
\end{multline*}

Note that by Lemma~\ref{LThetaCong}, we have that, for every $k\in\omega_{0}$, the following equality holds
\[
\mathrm{CH}^{(1)}_{s_{k}}\left(
\mathrm{ip}^{(1,X)@}_{s_{k}}\left(
R_{k}
\right)\right)
=
\mathrm{CH}^{(1)}_{s_{k}}\left(
\mathrm{ip}^{(1,X)@}_{s_{k}}\left(
R'_{k}
\right)\right)=
\mathrm{CH}^{(1)}_{s_{k}}\left(
\mathrm{ip}^{(1,X)@}_{s_{k}}\left(
R''_{k}
\right)\right).
\]

Therefore, the sequence
\[
\left(
\left(\mathrm{CH}^{(1)}_{s_{k}}\left(\mathrm{ip}^{(1,X)@}_{s_{k}}(R_{k})\right),s_{k}\right)
\right)_{k\in\omega_{0}}
\]
is a strictly decreasing $\omega_{0}$-chain in $(\coprod\mathrm{Pth}_{\boldsymbol{\mathcal{A}}}, \leq_{\mathbf{Pth}_{\boldsymbol{\mathcal{A}}}})$ contradicting the fact that, by Proposition~\ref{POrdArt}, $(\coprod\mathrm{Pth}_{\boldsymbol{\mathcal{A}}}, \leq_{\mathbf{Pth}_{\boldsymbol{\mathcal{A}}}})$ is an Artinian partially ordered set.

It follows that  $(\coprod[\mathrm{PT}_{\boldsymbol{\mathcal{A}}}], \leq_{[\mathbf{PT}_{\boldsymbol{\mathcal{A}}}]})$ is an Artinian partially ordered set.

This concludes the proof.
\end{proof}

We close this section by giving a characterisation of the order introduced in Definition~\ref{DPTQOrd}.

\begin{corollary}\label{CPTQOrd}
Let $s,t$ be sorts in $S$ and let $P,Q$ be path terms, with $P\in\mathrm{PT}_{\boldsymbol{\mathcal{A}},s}$ and $Q\in\mathrm{PT}_{\boldsymbol{\mathcal{A}},t}$. Then the following statements are equivalent.
\begin{enumerate}
\item $([Q]_{t},t)\leq_{[\mathbf{PT}_{\boldsymbol{\mathcal{A}}}]} ([P]_{s},s)$;
\item $(\mathrm{CH}^{(1)}_{t}(\mathrm{ip}^{(1,X)@}_{t}(Q)), t)
\leq_{\mathbf{T}_{\Sigma^{\boldsymbol{\mathcal{A}}}}(X)}
(\mathrm{CH}^{(1)}_{s}(\mathrm{ip}^{(1,X)@}_{s}(P)), s).
$
\end{enumerate}
\end{corollary}
\begin{proof}
Assume that the following inequality holds
\[
\left([Q]_{t},t\right)
\leq_{[\mathbf{PT}_{\boldsymbol{\mathcal{A}}}]} 
\left([P]_{s},s\right).
\]
Thus, by Definition~\ref{DPTQOrd}, there exists path term representatives $Q'\in [Q]_{t}$ and $P'\in [P]_{s}$ satisfying that 
\[
\left(
\mathrm{ip}^{(1,X)@}_{t}\left(
Q'
\right),
t
\right)
\leq_{\mathbf{Pth}_{\boldsymbol{\mathcal{A}}}}
\left(
\mathrm{ip}^{(1,X)@}_{s}\left(
P'
\right),
s
\right).
\]

Following Proposition~\ref{PCHMono}, we conclude that
\[
\left(
\mathrm{CH}^{(1)}_{t}\left(
\mathrm{ip}^{(1,X)@}_{t}\left(
Q'
\right)\right),t
\right)
\leq_{\mathbf{T}_{\Sigma^{\boldsymbol{\mathcal{A}}}}(X)}
\left(
\mathrm{CH}^{(1)}_{s}\left(
\mathrm{ip}^{(1,X)@}_{s}\left(
P'
\right)\right),s
\right).
\]

Note that, by Lemma~\ref{LThetaCong}, we have that
\allowdisplaybreaks
\begin{align*}
\mathrm{CH}^{(1)}_{t}\left(
\mathrm{ip}^{(1,X)@}_{t}\left(
Q
\right)\right)
&=
\mathrm{CH}^{(1)}_{t}\left(
\mathrm{ip}^{(1,X)@}_{t}\left(
Q'
\right)\right);
\\
\mathrm{CH}^{(1)}_{s}\left(
\mathrm{ip}^{(1,X)@}_{s}\left(
P
\right)\right)
&=
\mathrm{CH}^{(1)}_{s}\left(
\mathrm{ip}^{(1,X)@}_{s}\left(
P'
\right)\right).
\end{align*}

Now assume that the following inequality holds
\[
\left(
\mathrm{CH}^{(1)}_{t}\left(
\mathrm{ip}^{(1,X)@}_{t}\left(
Q
\right)\right), t\right)
\leq_{\mathbf{T}_{\Sigma^{\boldsymbol{\mathcal{A}}}}(X)}
\left(
\mathrm{CH}^{(1)}_{s}\left(
\mathrm{ip}^{(1,X)@}_{s}\left(
P
\right)\right), s\right).
\]

Then, by Proposition~\ref{PIpCHOrd}, the following inequality holds
\allowdisplaybreaks
\begin{multline*}
\left(
\mathrm{ip}^{(1,X)@}_{t}\left(
\mathrm{CH}^{(1)}_{t}\left(
\mathrm{ip}^{(1,X)@}_{t}\left(
Q
\right)\right)\right), t\right)
\\
\leq_{\mathbf{T}_{\Sigma^{\boldsymbol{\mathcal{A}}}}(X)}
\left(
\mathrm{ip}^{(1,X)@}_{s}\left(
\mathrm{CH}^{(1)}_{s}\left(
\mathrm{ip}^{(1,X)@}_{s}\left(
P
\right)\right)\right), s\right).
\end{multline*}

Note that, by Lemma~\ref{LWCong}, we have that 
\allowdisplaybreaks
\begin{align*}
\mathrm{CH}^{(1)}_{t}\left(
\mathrm{ip}^{(1,X)@}_{t}\left(
Q
\right)\right)&\in [Q]_{t};
&
\mathrm{CH}^{(1)}_{s}\left(
\mathrm{ip}^{(1,X)@}_{s}\left(
P
\right)\right)&\in [P]_{s}.
\end{align*}

Thus, by Definition~\ref{DPTQOrd}, we have that
\[
\left([Q]_{t},t\right)
\leq_{[\mathbf{PT}_{\boldsymbol{\mathcal{A}}}]} 
\left([P]_{s},s\right).
\]

This concludes the proof.
\end{proof}

\chapter{First-order isomorphisms}\label{S1J}

In this chapter we prove that the partial $\Sigma^{\boldsymbol{\mathcal{A}}}$-algebras of path classes $[\mathbf{Pth}_{\boldsymbol{\mathcal{A}}}]$ and path term classes $[\mathbf{PT}_{\boldsymbol{\mathcal{A}}}]$ are isomorphic. Indeed the mappings $\mathrm{ip}^{([1],X)@}$ and $\mathrm{CH}^{[1]}$ introduced in the previous chapter form a pair of inverse $\Sigma^{\boldsymbol{\mathcal{A}}}$-isomorphisms. Moreover, we show that the many-sorted categorial $\Sigma$-algebras $[\mathsf{Pth}_{\boldsymbol{\mathcal{A}}}]$ and $[\mathsf{PT}_{\boldsymbol{\mathcal{A}}}]$ are also isomorphic, since the mappings $\mathrm{ip}^{([1],X)@}$ and $\mathrm{CH}^{[1]}$ also form a pair of inverse functors, i.e., of categorial $\Sigma$-isomorphisms. Finally, we show that the mappings $\mathrm{ip}^{([1],X)@}$ and $\mathrm{CH}^{[1]}$ also form a pair of inverse order-preserving mappings which, in turn, entails that the Artinian partial orders $\leq_{[\mathbf{PT}_{\boldsymbol{\mathcal{A}}}]}$ and $\leq_{[\mathbf{Pth}_{\boldsymbol{\mathcal{A}}}]}$ are isomorphic.


\section{An algebraic isomorphism}

We begin by showing that the many-sorted partial $\Sigma^{\boldsymbol{\mathcal{A}}}$-algebras $[\mathbf{Pth}_{\boldsymbol{\mathcal{A}}}]$ of path classes, introduced in Proposition~\ref{PCHCatAlg}, and $[\mathbf{PT}_{\boldsymbol{\mathcal{A}}}]$ of path term classes, introduced in Proposition~\ref{PPTQCatAlg}, are isomorphic.

\begin{restatable}{theorem}{TIso}
\label{TIso} 
The partial $\Sigma^{\boldsymbol{\mathcal{A}}}$-algebras 
$[\mathbf{Pth}_{\boldsymbol{\mathcal{A}}}]$ and $[\mathbf{PT}_{\boldsymbol{\mathcal{A}}}]$ are isomorphic.
\end{restatable}
\begin{proof} We claim that the pair of $S$-sorted mappings
\begin{center}
\begin{tikzpicture}
[ACliment/.style={-{To [angle'=45, length=5.75pt, width=4pt, round]}}
, scale=0.8, 
AClimentD/.style={double equal sign distance,
-implies
}
]
\node[] (PK) at (-2.6,-3) [] {$[\mathbf{Pth}_{\boldsymbol{\mathcal{A}}}]$};
\node[] (PTK) at (2.6,-3) [] {$[\mathbf{PT}_{\boldsymbol{\mathcal{A}}}]$};
 
 \draw[shorten >=0.2cm,shorten <=.2cm, ACliment]
 ($(-2,-3)+.09*(0,3)$) to node  [above]  {$\mathrm{CH}^{[1]}$}  ($(2,-3)+.09*(0,3)$);
 \draw[shorten >=0.2cm,shorten <=.2cm, ACliment]
 ($(2,-3)-.09*(0,3)$) to node  [below]  {$\mathrm{ip}^{([1],X)@}$}  ($(-2,-3)-.09*(0,3)$);
\end{tikzpicture}
\end{center}
determines a pair of $\Sigma^{\boldsymbol{\mathcal{A}}}$-isomorphisms such that 
\allowdisplaybreaks
\begin{align*}
\mathrm{ip}^{([1],X)@}
\circ
\mathrm{CH}^{[1]}
&=
\mathrm{id}^{[\mathbf{Pth}_{\boldsymbol{\mathcal{A}}}]};
&
\mathrm{CH}^{[1]}
\circ
\mathrm{ip}^{([1],X)@}
&=
\mathrm{id}^{[\mathbf{PT}_{\boldsymbol{\mathcal{A}}}]}.
\end{align*}

We prove some intermediate claims that will allow us to conclude the aforementioned statement.

\begin{restatable}{claim}{CIsoCH}
\label{CIsoCH}
The mapping
$\mathrm{CH}^{[1]}$
is a $\Sigma^{\boldsymbol{\mathcal{A}}}$-homomorphism.
\end{restatable}

Let us recall that $\mathrm{CH}^{[1]}$, introduced in Definition~\ref{DPTQCH}, is a many-sorted mapping of the form
$$
\mathrm{CH}^{[1]}
\colon
[\mathrm{Pth}_{\boldsymbol{\mathcal{A}}}]
\mor
[\mathrm{PT}_{\boldsymbol{\mathcal{A}}}].
$$

By Definition~\ref{DPTQCH}, the many-sorted mapping $\mathrm{CH}^{[1]}$ is given by
\[
\mathrm{CH}^{[1]} = \mathrm{pr}^{\Theta^{[1]}}\circ \mathrm{CH}^{(1)\mathrm{m}} = \mathrm{pr}^{\Theta^{[1]}}\circ \mathrm{CH}^{(1)}\circ \mathrm{pr}^{\mathrm{Ker}(\mathrm{CH}^{(1)})}.
\]

By Proposition~\ref{PCHCatAlg}, $\mathrm{pr}^{\mathrm{Ker}(\mathrm{CH}^{(1)})}$ is a $\Sigma^{\boldsymbol{\mathcal{A}}}$-homomorphism. Following Proposition~\ref{PThetaCH}, $\mathrm{pr}^{\Theta^{[1]}}\circ \mathrm{CH}^{(1)}$ is a $\Sigma^{\boldsymbol{\mathcal{A}}}$-homomorphism. It follows that $\mathrm{CH}^{[1]}$, being a composition of $\Sigma^{\boldsymbol{\mathcal{A}}}$-homomorphisms, is itself a $\Sigma^{\boldsymbol{\mathcal{A}}}$-homomorphism.

This finishes the proof of Claim~\ref{CIsoCH}.

Now, we turn ourselves to study the case of the many-sorted mapping $\mathrm{ip}^{([1],X)@}$.

\begin{restatable}{claim}{CIsoIpfc}
\label{CIsoIpfc} 
The mapping $\mathrm{ip}^{([1],X)@}$ is a $\Sigma^{\boldsymbol{\mathcal{A}}}$-homomorphism.
\end{restatable}

Let us recall that $\mathrm{ip}^{([1],X)@}$, introduced in Definition~\ref{DPTQIp}, is a many-sorted mapping of the form
$$
\mathrm{ip}^{([1],X)@}
\colon
\left[\mathrm{PT}_{\boldsymbol{\mathcal{A}}}\right]
\mor
\left[ \mathrm{Pth}_{\boldsymbol{\mathcal{A}}}\right].
$$

By Definition~\ref{DPTQIp}, the many-sorted mapping $
\mathrm{ip}^{([1],X)@}$ is given by
\[
(\mathrm{pr}^{\mathrm{Ker}(\mathrm{CH}^{(1)})}\circ\mathrm{ip}^{(1,X)@})^{\mathrm{m}}\circ\mathrm{p}^{ \mathrm{Ker}(\mathrm{pr}^{\mathrm{Ker}(\mathrm{CH}^{(1)})}\circ\mathrm{ip}^{(1,X)@}),\Theta^{[1]}}.
\]

By Definition~\ref{DPTQIp}, this mapping satisfies that 
\[
\mathrm{pr}^{\mathrm{Ker}(\mathrm{CH}^{(1)})}\circ\mathrm{ip}^{(1,X)@}
=
\mathrm{ip}^{([1],X)@}
\circ
\mathrm{pr}^{\Theta^{[1]}}.
\]

Note that $\mathrm{ip}^{(1,X)@}$ is a $\Sigma^{\boldsymbol{\mathcal{A}}}$-homomorphism according to Definition~\ref{DIp} and $\mathrm{pr}^{\mathrm{Ker}(\mathrm{CH}^{(1)})}$ is a $\Sigma^{\boldsymbol{\mathcal{A}}}$-homomorphism according to Proposition~\ref{PCHCatAlg}. It follows that  $\mathrm{ip}^{([1],X)@}$, being a composition of $\Sigma^{\boldsymbol{\mathcal{A}}}$-homomorphisms, is itself a $\Sigma^{\boldsymbol{\mathcal{A}}}$-homomorphism.

This finishes the proof of Claim~\ref{CIsoIp}.

\begin{restatable}{claim}{CIso}
 \label{CIso}
$\mathrm{CH}^{[1]}$ and $\mathrm{ip}^{([1],X)@}$ form a pair of inverse $\Sigma^{\boldsymbol{\mathcal{A}}}$-isomorphisms.
\end{restatable}

Let $s\in S$ and $[\mathfrak{P}]_{s}$ be a path class in $[\mathrm{Pth}_{\boldsymbol{\mathcal{A}}}]_{s}$.

The following chain of equalities holds
\begin{align*}
\mathrm{ip}^{([1],X)@}_{s}\left(
\mathrm{CH}^{[1]}_{s}\left(
\left[
\mathfrak{P}
\right]_{s}
\right)\right)
&=
\left[
\mathrm{ip}^{(1,X)@}_{s}\left(
\mathrm{CH}^{(1)}_{s}\left(
\mathfrak{P}
\right)\right)\right]_{s}
\tag{1}
\\&=
\left[\mathfrak{P}\right]_{s}.
\tag{2}
\end{align*}

The first equality unravels the definition of the mappings $\mathrm{ip}^{([1],X)@}$ and $\mathrm{CH}^{[1]}$ according to, respectively, Definitions~\ref{DPTQIp} and~\ref{DPTQCH}; finally, the second equality follows from Proposition~\ref{PIpCH}.

On the other hand, let $s$ be a sort in $S$ and $[P]_{s}$  a path term class in $[\mathrm{PT}_{\boldsymbol{\mathcal{A}}}]_{s}$.

The following chain of equalities holds
\allowdisplaybreaks
\begin{align*}
\mathrm{CH}^{[1]}_{s}\left(
\mathrm{ip}^{([1],X)@}_{s}\left(
\left[
P
\right]_{s}
\right)\right)
&=
\left[
\mathrm{CH}^{(1)}_{s}\left(
\mathrm{ip}^{(1,X)@}_{s}\left(
P
\right)\right)\right]_{s}
\tag{1}
\\&=
\left[
P
\right]_{s}.
\tag{2}
\end{align*}

The first equality unravels the definition of the mappings $\mathrm{ip}^{([1],X)@}$ and $\mathrm{CH}^{[1]}$ according to, respectively, Definitions~\ref{DPTQIp} and~\ref{DPTQCH}; finally, the last equation follows from Lemma~\ref{LWCong}.

This completes the proof of Claim~\ref{CIso}.

This completes the proof of Theorem~\ref{TIso}.
\end{proof}

\section{A categorial isomorphism}

In this section we prove that the $S$-sorted categories $[\mathsf{Pth}_{\boldsymbol{\mathcal{A}}}]$, introduced in Proposition~\ref{PCHCat}, and $[\mathsf{PT}_{\boldsymbol{\mathcal{A}}}]$, introduced in Definition~\ref{DPTQCat}, are isomorphic. 

\begin{restatable}{theorem}{TIsoCat}\label{TIsoCat} The $S$-sorted categories $[\mathsf{Pth}_{\boldsymbol{\mathcal{A}}}]$ and $[\mathsf{PT}_{\boldsymbol{\mathcal{A}}}]$ are isomorphic.
\end{restatable}
\begin{proof}
Consider the $S$-sorted mapping $\mathrm{CH}^{[1]}$ introduced in Definition~\ref{DPTQCH}, where
\[
\mathrm{CH}^{[1]}\colon [\mathbf{Pth}_{\boldsymbol{\mathcal{A}}}]\mor [\mathbf{PT}_{\boldsymbol{\mathcal{A}}}]
\]

\begin{restatable}{claim}{CIsoCatCH}
\label{CIsoCatCH} $\mathrm{CH}^{[1]}$ is a functor from $[\mathsf{Pth}_{\boldsymbol{\mathcal{A}}}]$ to $[\mathsf{PT}_{\boldsymbol{\mathcal{A}}}]$.
\end{restatable}
Following Definition~\ref{DCat}, we need to prove that the following equalities hold

\textsf{(1)} Let $s$ be a sort in $S$ and $[\mathfrak{P}]_{s}$ a path class in $[\mathrm{Pth}_{\boldsymbol{\mathcal{A}}}]_{s}$, then the following equalities hold
\allowdisplaybreaks
\begin{align*}
\mathrm{CH}^{[1]}_{s}\left(
\mathrm{sc}^{0[\mathbf{Pth}_{\boldsymbol{\mathcal{A}}}]}_{s}\left(
[\mathfrak{P}]_{s}
\right)
\right)
&=
\mathrm{sc}^{0[\mathbf{PT}_{\boldsymbol{\mathcal{A}}}]}_{s}\left(
\mathrm{CH}^{[1]}_{s}\left(
[\mathfrak{P}]_{s}
\right)
\right);
\\
\mathrm{CH}^{[1]}_{s}\left(
\mathrm{tg}^{0[\mathbf{Pth}_{\boldsymbol{\mathcal{A}}}]}_{s}\left(
[\mathfrak{P}]_{s}
\right)
\right)
&=
\mathrm{tg}^{0[\mathbf{PT}_{\boldsymbol{\mathcal{A}}}]}_{s}\left(
\mathrm{CH}^{[1]}_{s}\left(
[\mathfrak{P}]_{s}
\right)
\right).
\end{align*}

\textsf{(2)} Let $s$ be a sort in $S$ and $[\mathfrak{Q}]_{s}, [\mathfrak{P}]_{s}$ be path classes in $[\mathrm{Pth}_{\boldsymbol{\mathcal{A}}}]_{s}$ satisfying that 
\[
\mathrm{sc}^{0[\mathbf{Pth}_{\boldsymbol{\mathcal{A}}}]}_{s}\left([\mathfrak{Q}]_{s}\right)
=
\mathrm{tg}^{0[\mathbf{Pth}_{\boldsymbol{\mathcal{A}}}]}_{s}\left([\mathfrak{P}]_{s}\right),
\]
then the following equality holds
\[
\mathrm{CH}^{[1]}_{s}\left(
[\mathfrak{Q}]_{s}
\circ^{0[\mathbf{Pth}_{\boldsymbol{\mathcal{A}}}]}_{s}
[\mathfrak{P}]_{s}
\right)=
\mathrm{CH}^{[1]}_{s}\left(
[\mathfrak{Q}]_{s}
\right)
\circ^{0[\mathbf{PT}_{\boldsymbol{\mathcal{A}}}]}_{s}
\mathrm{CH}^{[1]}_{s}\left(
[\mathfrak{P}]_{s}
\right).
\]

All these equalities follow from the fact that $\mathrm{CH}^{[1]}$ is a many-sorted partial $\Sigma^{\boldsymbol{\mathcal{A}}}$-homomorphism from $[\mathbf{Pth}_{\boldsymbol{\mathcal{A}}}]$ to $[\mathbf{PT}_{\boldsymbol{\mathcal{A}}}]$ according to Theorem~\ref{TIso}.

This completes the proof of Claim~\ref{CIsoCatCH}.

Now consider the $S$-sorted mapping $\mathrm{ip}^{([1],X)@}$ introduced in Definition~\ref{DPTQIp}, where
\[
\mathrm{ip}^{([1],X)@}\colon [\mathbf{PT}_{\boldsymbol{\mathcal{A}}}]
\mor [\mathbf{Pth}_{\boldsymbol{\mathcal{A}}}]
\]

\begin{restatable}{claim}{CIsoCatIp}
\label{CIsoCatIp} $\mathrm{ip}^{([1],X)@}$ is a functor from $[\mathsf{PT}_{\boldsymbol{\mathcal{A}}}]$ to $[\mathsf{Pth}_{\boldsymbol{\mathcal{A}}}]$.
\end{restatable}
Following Definition~\ref{DCat}, we need to prove that the following equalities hold

\textsf{(1)} Let $s$ be a sort in $S$ and $[P]_{s}$ a path term class in $[\mathrm{PT}_{\boldsymbol{\mathcal{A}}}]_{s}$, then the following equalities hold
\allowdisplaybreaks
\begin{align*}
\mathrm{ip}^{([1],X)@}_{s}\left(
\mathrm{sc}^{0[\mathbf{PT}_{\boldsymbol{\mathcal{A}}}]}_{s}\left(
[P]_{s}
\right)
\right)
&=
\mathrm{sc}^{0[\mathbf{Pth}_{\boldsymbol{\mathcal{A}}}]}_{s}\left(
\mathrm{ip}^{([1],X)@}_{s}\left(
[P]_{s}
\right)
\right);
\\
\mathrm{ip}^{([1],X)@}_{s}\left(
\mathrm{tg}^{0[\mathbf{PT}_{\boldsymbol{\mathcal{A}}}]}_{s}\left(
[P]_{s}
\right)
\right)
&=
\mathrm{tg}^{0[\mathbf{Pth}_{\boldsymbol{\mathcal{A}}}]}_{s}\left(
\mathrm{ip}^{([1],X)@}_{s}\left(
[P]_{s}
\right)
\right).
\end{align*}

\textsf{(2)} Let $s$ be a sort in $S$ and $[Q]_{s}, [P]_{s}$ be path term classes in $[\mathrm{PT}_{\boldsymbol{\mathcal{A}}}]_{s}$ satisfying that 
\[
\mathrm{sc}^{0[\mathbf{PT}_{\boldsymbol{\mathcal{A}}}]}_{s}\left(
[Q]_{s}\right)
=
\mathrm{tg}^{0[\mathbf{PT}_{\boldsymbol{\mathcal{A}}}]}_{s}\left(
[P]_{s}\right),
\]
then the following equality holds
\[
\mathrm{ip}^{([1],X)@}_{s}\left(
[Q]_{s}
\circ^{0[\mathbf{PT}_{\boldsymbol{\mathcal{A}}}]}_{s}
[P]_{s}
\right)=
\mathrm{ip}^{([1],X)@}_{s}\left(
[Q]_{s}
\right)
\circ^{0[\mathbf{Pth}_{\boldsymbol{\mathcal{A}}}]}_{s}
\mathrm{ip}^{([1],X)@}_{s}\left(
[P]_{s}
\right).
\]

All these equalities follow from the fact that $\mathrm{ip}^{([1],X)@}$ is a many-sorted partial $\Sigma^{\boldsymbol{\mathcal{A}}}$-homomorphism from $[\mathbf{PT}_{\boldsymbol{\mathcal{A}}}]$ to $[\mathbf{Pth}_{\boldsymbol{\mathcal{A}}}]$ according to Theorem~\ref{TIso}.

This completes the proof of Claim~\ref{CIsoCatIp}.

That the following equalities between functors hold is a direct consequence of Theorem~\ref{TIso}.
\begin{align*}
\mathrm{ip}^{([1],X)@}\circ \mathrm{CH}^{[1]}&=
\mathrm{id}^{[\mathsf{Pth}_{\boldsymbol{\mathcal{A}}}]};
&
\mathrm{CH}^{[1]}\circ \mathrm{ip}^{([1],X)@}&=
\mathrm{id}^{[\mathsf{PT}_{\boldsymbol{\mathcal{A}}}]}.
\end{align*}

This concludes the proof.
\end{proof}

We next prove the compatibility of the above functors with the underlying $\Sigma$-algebraic structures of  $[\mathsf{Pth}_{\boldsymbol{\mathcal{A}}}]$ and $[\mathsf{PT}_{\boldsymbol{\mathcal{A}}}]$.

\begin{restatable}{lemma}{CIsoCatAlgCH}
\label{CIsoCatAlgCH} $\mathrm{CH}^{[1]}$ is a categorial $\Sigma$-homomorphism from $[\mathsf{Pth}_{\boldsymbol{\mathcal{A}}}]$ to $[\mathsf{PT}_{\boldsymbol{\mathcal{A}}}]$.
\end{restatable}
\begin{proof}
By Definition~\ref{DnCatAlg}, we need to prove that, for every pair $(\mathbf{s},s)$ in $\Sigma_{\mathbf{s},s}$, every operation $\sigma\in \Sigma_{\mathbf{s},s}$ and every family of path classes $([\mathfrak{P}_{j}]_{s_{j}})_{j\in\bb{\mathbf{s}}}$ in $[\mathrm{Pth}_{\boldsymbol{\mathcal{A}}}]_{\mathbf{s}}$, the following equality holds
\[
\mathrm{CH}^{[1]}_{s}\left(
\sigma^{[\mathbf{Pth}_{\boldsymbol{\mathcal{A}}}]}\left(
\left(
[\mathfrak{P}_{j}]_{s_{j}}
\right)_{j\in\bb{\mathbf{s}}}
\right)
\right)
=
\sigma^{[\mathbf{PT}_{\boldsymbol{\mathcal{A}}}]}\left(
\left(
\mathrm{CH}^{[1]}_{s_{j}}\left(
[\mathfrak{P}_{j}]_{s_{j}}
\right)
\right)_{j\in\bb{\mathbf{s}}}
\right).
\]

This equality follows from the fact that $\mathrm{CH}^{[1]}$ is a many-sorted partial $\Sigma^{\boldsymbol{\mathcal{A}}}$-homomorphism from $[\mathbf{Pth}_{\boldsymbol{\mathcal{A}}}]$ to $[\mathbf{PT}_{\boldsymbol{\mathcal{A}}}]$ according to Theorem~\ref{TIso}.

This completes the proof.
\end{proof}

\begin{restatable}{lemma}{CIsoCatAlgIp}
\label{CIsoCatAlgIp} $\mathrm{ip}^{([1],X)@}$ is a categorial $\Sigma$-homomorphism from  $[\mathsf{PT}_{\boldsymbol{\mathcal{A}}}]$ to $[\mathsf{Pth}_{\boldsymbol{\mathcal{A}}}]$.
\end{restatable}

\begin{proof}
By Definition~\ref{DnCatAlg}, we need to prove that, for every pair $(\mathbf{s},s)$ in $\Sigma_{\mathbf{s},s}$, every operation $\sigma\in \Sigma_{\mathbf{s},s}$ and every family of path term classes $([P_{j}]_{s_{j}})_{j\in\bb{\mathbf{s}}}$ in $[\mathrm{PT}_{\boldsymbol{\mathcal{A}}}]_{\mathbf{s}}$, the following equality holds
\[
\mathrm{ip}^{([1],X)@}_{s}\left(
\sigma^{[\mathbf{PT}_{\boldsymbol{\mathcal{A}}}]}\left(
\left(
[P_{j}]_{s_{j}}
\right)_{j\in\bb{\mathbf{s}}}
\right)
\right)
=
\sigma^{[\mathbf{Pth}_{\boldsymbol{\mathcal{A}}}]}\left(
\left(
\mathrm{ip}^{([1],X)@}_{s_{j}}\left(
[P_{j}]_{s_{j}}
\right)
\right)_{j\in\bb{\mathbf{s}}}
\right).
\]

This equality follows from the fact that $\mathrm{ip}^{([1],X)@}$ is a many-sorted partial $\Sigma^{\boldsymbol{\mathcal{A}}}$-homomorphism from  $[\mathbf{PT}_{\boldsymbol{\mathcal{A}}}]$ to $[\mathbf{Pth}_{\boldsymbol{\mathcal{A}}}]$ according to Theorem~\ref{TIso}.

This completes the proof.
\end{proof}

As a consequence of the above results we can conclude that the $S$-sorted categorial $\Sigma$-algebras of path classes $[\mathsf{Pth}_{\boldsymbol{\mathcal{A}}}]$ and of path term classes $[\mathsf{PT}_{\boldsymbol{\mathcal{A}}}]$ are isomorphic. 

\begin{restatable}{corollary}{CIsoCatAlg}
\label{CIsoCatAlg} The $S$-sorted categorial $\Sigma$-algebras $[\mathsf{Pth}_{\boldsymbol{\mathcal{A}}}]$ and $[\mathsf{PT}_{\boldsymbol{\mathcal{A}}}]$ are isomorphic.
\end{restatable}
\begin{proof}
From Theorem~\ref{TIsoCat} we have that the underlying $S$-sorted categories of $[\mathsf{Pth}_{\boldsymbol{\mathcal{A}}}]$ and $[\mathsf{PT}_{\boldsymbol{\mathcal{A}}}]$ are isomorphic by means of the pair of functors
\begin{center}
\begin{tikzpicture}
[ACliment/.style={-{To [angle'=45, length=5.75pt, width=4pt, round]}}
, scale=0.8, 
AClimentD/.style={double equal sign distance,
-implies
}
]
\node[] (PK) at (-2.8,-3) [] {$[\mathsf{Pth}_{\boldsymbol{\mathcal{A}}}]$};
\node[] (PTK) at (2.8,-3) [] {$[\mathsf{PT}_{\boldsymbol{\mathcal{A}}}]$};
 
 \draw[shorten >=0.2cm,shorten <=.2cm, ACliment]
 ($(-2,-3)+.09*(0,3)$) to node  [above]  {$\mathrm{CH}^{[1]}$}  ($(2,-3)+.09*(0,3)$);
 \draw[shorten >=0.2cm,shorten <=.2cm, ACliment]
 ($(2,-3)-.09*(0,3)$) to node  [below]  {$\mathrm{ip}^{([1],X)@}$}  ($(-2,-3)-.09*(0,3)$);
\end{tikzpicture}
\end{center}

Moreover, by Lemma~\ref{CIsoCatAlgCH}, $\mathrm{CH}^{[1]}$ is a categorial $\Sigma$-homomorphism from $[\mathsf{Pth}_{\boldsymbol{\mathcal{A}}}]$ to $[\mathsf{PT}_{\boldsymbol{\mathcal{A}}}]$, and, by 
Lemma~\ref{CIsoCatAlgIp}, $\mathrm{ip}^{([1],X)@}$ is a categorial $\Sigma$-homomorphism from  $[\mathsf{PT}_{\boldsymbol{\mathcal{A}}}]$ to $[\mathsf{Pth}_{\boldsymbol{\mathcal{A}}}]$.

This concludes the proof.
\end{proof}

\section{An order isomorphism}

In this section we prove that the coproduct of $\mathrm{CH}^{[1]}$ and the coproduct of $\mathrm{ip}^{([1],X)@}$ are order-preserving mappings. This will ultimately lead to prove that the Artinian partially ordered sets defined, on one hand, on Definition~\ref{DCHOrd}, on the coproduct of  path classes and, on the other hand, on Definition~\ref{DPTQOrd}, on the coproduct of path term classes  are order isomorphic.

We start by proving that the coproduct of the mapping $\mathrm{CH}^{[1]}$ is order-preserving.

\begin{restatable}{lemma}{LIsoOrdCH}
\label{LIsoOrdCH} 
The mapping $\coprod\mathrm{CH}^{[1]}$ from $\coprod[\mathrm{Pth}_{\boldsymbol{\mathcal{A}}}]$ to  $\coprod[\mathrm{PT}_{\boldsymbol{\mathcal{A}}}]$ that, for every sort $s$ in $S$ and every path class $[\mathfrak{P}]_{s}$ in $[\mathrm{Pth}_{\boldsymbol{\mathcal{A}}}]_{s}$ sends the pair $([\mathfrak{P}]_{s},s)$ in $\coprod[\mathrm{Pth}_{\boldsymbol{\mathcal{A}}}]$ to $([\mathrm{CH}^{(1)}_{s}(\mathfrak{P})]_{s},s)$ in $\coprod[\mathrm{PT}_{\boldsymbol{\mathcal{A}}}]$ determines an order-preserving mapping
\[
\textstyle
\coprod\mathrm{CH}^{[1]}\colon
\Big(\coprod[\mathrm{Pth}_{\boldsymbol{\mathcal{A}}}]
,\,\leq_{[\mathbf{Pth}_{\boldsymbol{\mathcal{A}}}]}
\Big)
\mor
\textstyle
\Big(\coprod[\mathrm{PT}_{\boldsymbol{\mathcal{A}}}],\,
\leq_{[\mathbf{PT}_{\boldsymbol{\mathcal{A}}}]}
\Big).
\]
\end{restatable}
\begin{proof}
Let $s,t$ be sorts in $S$ and let us consider paths $\mathfrak{Q}$ in $\mathrm{Pth}_{\boldsymbol{\mathcal{A}},t}$ and $\mathfrak{P}$ in $\mathrm{Pth}_{\boldsymbol{\mathcal{A}},s}$ satisfying that
$
\left(\left[
\mathfrak{Q}
\right]_{t},t\right)
\leq_{[\mathbf{Pth}_{\boldsymbol{\mathcal{A}}}]}
\left(\left[\mathfrak{P}
\right]_{s},s\right)
$. 
We have to prove that 
$
\textstyle
\coprod\mathrm{CH}^{[1]}([\mathfrak{Q}]_{t},t)
\leq_{[\mathbf{PT}_{\boldsymbol{\mathcal{A}}}]}
\coprod\mathrm{CH}^{[1]}_{s}([\mathfrak{P}]_{s},s)
$. 
Taking into account the definition of the mapping $\coprod\mathrm{CH}^{[1]}$, introduced in Definition~\ref{DPTQCH},  the presentation of the partial order $\leq_{[\mathbf{PT}_{\boldsymbol{\mathcal{A}}}]}$, introduced in Proposition~\ref{DPTQOrd}, Corollary~\ref{CPTQOrd} and Proposition~\ref{PIpCH}, this is equivalent to prove the following inequality
$$
\left(
\mathrm{CH}^{(1)}_{t}\left(
\mathfrak{Q}
\right),t\right)
\leq_{\mathbf{T}_{\Sigma^{\boldsymbol{\mathcal{A}}}}(X)}
\left(\mathrm{CH}^{(1)}_{s}\left(
\mathfrak{P}\right),s\right).
$$

Now, since $
([\mathfrak{Q}]_{t},t)\leq_{[\mathbf{Pth}_{\boldsymbol{\mathcal{A}}}]} ([\mathfrak{P}]_{s},s)$ we have, according to Definition~\ref{DCHOrd}, that there exists a path $\mathfrak{Q}'$ in $[\mathfrak{Q}]_{t}$ and a path $\mathfrak{P}'\in[\mathfrak{P}]_{s}$ for which
$$
\left(
\mathfrak{Q}',t
\right)
\leq_{\mathbf{Pth}_{\boldsymbol{\mathcal{A}}}}
\left(\mathfrak{P}',s\right).
$$
 
Note that $\mathrm{CH}^{(1)}_{t}(\mathfrak{Q}')=\mathrm{CH}^{(1)}_{t}(\mathfrak{Q})$ and $\mathrm{CH}^{(1)}_{s}(\mathfrak{P}')=\mathrm{CH}^{(1)}_{s}(\mathfrak{P})$.

Since  
$
(\mathfrak{Q}',t)
\leq_{\mathbf{Pth}_{\boldsymbol{\mathcal{A}}}}
 (\mathfrak{P}',s)
$ 
and taking into account, by Proposition~\ref{PCHMono}, that $\coprod\mathrm{CH}^{(1)}$ is an order-preserving mapping, we have that 
$$
\left(
\mathrm{CH}^{(1)}_{t}\left(
\mathfrak{Q}'
\right),t\right)
\leq_{\mathbf{T}_{\Sigma^{\boldsymbol{\mathcal{A}}}}(X)}
\left(
\mathrm{CH}^{(1)}_{s}\left(
\mathfrak{P}'
\right),s
\right).
$$

This completes the proof.
\end{proof}

We next prove that the coproduct of the mapping $\mathrm{ip}^{([1],X)@}$ is order-preserving.

\begin{restatable}{lemma}{LIsoOrdIp}
\label{LIsoOrdIp} 
The mapping $\coprod\mathrm{ip}^{([1],X)@}$ from $\coprod[\mathrm{PT}_{\boldsymbol{\mathcal{A}}}]$ to  $\coprod[\mathrm{Pth}_{\boldsymbol{\mathcal{A}}}]$  that, for every sort $s$ in $S$ and every path term class $[P]_{s}$ in $[\mathrm{PT}_{\boldsymbol{\mathcal{A}}}]_{s}$, sends $([P]_{s},s)$ in $\coprod[\mathrm{PT}_{\boldsymbol{\mathcal{A}}}]_{s}$ to $([\mathrm{ip}^{(1,X)@}_{s}(P)]_{s},s)$ in $\coprod [\mathrm{Pth}_{\boldsymbol{\mathcal{A}}}]$ determines an order-preserving mapping
\[
\textstyle
\coprod\mathrm{ip}^{([1],X)@}\colon
\Big(\coprod[\mathrm{PT}_{\boldsymbol{\mathcal{A}}}],\,
\leq_{[\mathbf{PT}_{\boldsymbol{\mathcal{A}}}]}
\Big)
\mor
\textstyle
\Big(\coprod[\mathrm{Pth}_{\boldsymbol{\mathcal{A}}}]
,\,\leq_{[\mathbf{Pth}_{\boldsymbol{\mathcal{A}}}]}\Big).
\]
\end{restatable}
\begin{proof}
Let $s,t$ be sorts in $S$ and let us consider path terms $Q$ in $\mathrm{PT}_{\boldsymbol{\mathcal{A}}, t}$ and $P$ in $\mathrm{PT}_{\boldsymbol{\mathcal{A}}, s}$ satisfying that 
$
\left(\left[
Q
\right]_{t},t\right)
\leq_{[\mathbf{PT}_{\boldsymbol{\mathcal{A}}}]}
\left(\left[
P
\right]_{s},s\right)
$. 
We have to prove that $
\textstyle
\coprod\mathrm{ip}^{([1],X)@}([Q]_{t},t)
\leq_{[\mathbf{Pth}_{\boldsymbol{\mathcal{A}}}]}
\coprod\mathrm{ip}^{([1],X)@}([P]_{s},s)
$. Taking into account the definition of the mapping $\coprod\mathrm{ip}^{([1],X)@}$ and Definition~\ref{DCHOrd} this is equivalent to prove the existence of a path $\mathfrak{Q}'\in [\mathrm{ip}^{(1,X)@}_{t}(Q)]_{t}$ and a path $\mathfrak{P}'\in[\mathrm{ip}^{(1,X)@}_{s}(P)]_{s}$ satisfying the  inequality
$
\left(\mathfrak{Q}',t\right)
\leq_{\mathbf{Pth}_{\boldsymbol{\mathcal{A}}}}
\left(\mathfrak{P}',s\right)
$.

Note that, since 
$
([Q]_{t},t)
\leq_{[\mathbf{PT}_{\boldsymbol{\mathcal{A}}}]}
([P]_{s},s)$, we have, according to Corollary~\ref{CPTQOrd}, that
$$
\left(
\mathrm{CH}^{(1)}_{t}\left(
\mathrm{ip}^{(1,X)@}_{t}\left(
Q
\right)\right),t\right)
\leq_{\mathbf{T}_{\Sigma^{\boldsymbol{\mathcal{A}}}}(X)}
\left(
\mathrm{CH}^{(1)}_{s}\left(
\mathrm{ip}^{(1,X)@}_{s}\left(
P
\right)\right),s\right).
$$

Following Proposition~\ref{PIpCHOrd}, we have that
\allowdisplaybreaks
\begin{multline*}
\left(
\mathrm{ip}^{(1,X)@}_{t}\left(
\mathrm{CH}^{(1)}_{t}\left(
\mathrm{ip}^{(1,X)@}_{t}\left(
Q
\right)\right)\right),t\right)
\\
\leq_{\mathbf{Pth}_{\boldsymbol{\mathcal{A}}}}
\left(
\mathrm{ip}^{(1,X)@}_{s}\left(
\mathrm{CH}^{(1)}_{s}\left(
\mathrm{ip}^{(1,X)@}_{s}\left(
P
\right)\right)\right),s\right).
\end{multline*}

The statement follows from Lemmas~\ref{LWCong} and~\ref{LThetaCong}. Note that 
\allowdisplaybreaks
\begin{align*}
\mathrm{ip}^{(1,X)@}_{t}\left(
\mathrm{CH}^{(1)}_{t}\left(
\mathrm{ip}^{(1,X)@}_{t}\left(
Q
\right)\right)\right)&\in \left[\mathrm{ip}^{(1,X)@}_{t}\left(
Q
\right)\right]_{t};
\\
\mathrm{ip}^{(1,X)@}_{s}\left(
\mathrm{CH}^{(1)}_{s}\left(
\mathrm{ip}^{(1,X)@}_{s}\left(
P
\right)\right)\right)&\in \left[\mathrm{ip}^{(1,X)@}_{s}\left(
P
\right)\right]_{s}.
\end{align*}

This completes the proof.
\end{proof}

\begin{restatable}{theorem}{TIsoOrd}
\label{TIsoOrd} The Artinian ordered sets $(\coprod[\mathrm{Pth}_{\boldsymbol{\mathcal{A}}}]
,\,\leq_{[\mathbf{Pth}_{\boldsymbol{\mathcal{A}}}]})$, of path classes, and
$(\coprod[\mathrm{PT}_{\boldsymbol{\mathcal{A}}}],\leq_{[\mathbf{PT}_{\boldsymbol{\mathcal{A}}}]})$, of path term classes, are order isomorphic.
\end{restatable}
\begin{proof}
By Lemmas~\ref{LIsoOrdCH} and~\ref{LIsoOrdIp}, the mappings $\coprod\mathrm{CH}^{[1]}$ and $\coprod\mathrm{ip}^{([1],X)@}$ form a pair of mutually inverse order-preserving bijections.
\end{proof}

\chapter{
\texorpdfstring
{Freedom in the $\mathrm{QE}$-variety $\mathcal{V}(\boldsymbol{\mathcal{E}}^{\boldsymbol{\mathcal{A}}})$}
{Freedom in a variety}
}\label{S1K}

In this chapter we define, for the rewriting system $\boldsymbol{\mathcal{A}} = ((\mathbf{\Sigma},X),\mathcal{A})$, a specification $\boldsymbol{\mathcal{E}}^{\boldsymbol{\mathcal{A}}}$, whose defining equations 
$\mathcal{E}^{\boldsymbol{\mathcal{A}}}$ are $\mathrm{QE}$-equations, the $\mathrm{QE}$-variety of partial $\Sigma^{\boldsymbol{\mathcal{A}}}$-algebras $\mathcal{V}(\boldsymbol{\mathcal{E}}^{\boldsymbol{\mathcal{A}}})$ determined by $\boldsymbol{\mathcal{E}}^{\boldsymbol{\mathcal{A}}}$ and the category $\mathbf{PAlg}(\boldsymbol{\mathcal{E}}^{\boldsymbol{\mathcal{A}}})$ associated to $\mathcal{V}(\boldsymbol{\mathcal{E}}^{\boldsymbol{\mathcal{A}}})$. We show that the partial $\Sigma^{\boldsymbol{\mathcal{A}}}$-algebra of path classes $[\mathbf{Pth}_{\boldsymbol{\mathcal{A}}}]$ satisfies the defining equations $\mathcal{E}^{\boldsymbol{\mathcal{A}}}$ of the variety $\mathcal{V}(\boldsymbol{\mathcal{E}}^{\boldsymbol{\mathcal{A}}})$, i.e., that $[\mathbf{Pth}_{\boldsymbol{\mathcal{A}}}]\in \mathcal{V}(\boldsymbol{\mathcal{E}}^{\boldsymbol{\mathcal{A}}})$.  Then we prove the fundamental result of this chapter: that the two partial $\Sigma^{\boldsymbol{\mathcal{A}}}$-algebras $\mathbf{T}_{\boldsymbol{\mathcal{E}}^{\boldsymbol{\mathcal{A}}}}(\mathbf{Pth}_{\boldsymbol{\mathcal{A}}})$, which is the free partial $\Sigma^{\boldsymbol{\mathcal{A}}}$-algebra in the category $\mathbf{PAlg}(\boldsymbol{\mathcal{E}}^{\boldsymbol{\mathcal{A}}})$, and $[\mathbf{Pth}_{\boldsymbol{\mathcal{A}}}]$ are isomorphic. From here, as a consequence of the isomorphism of the previous chapter, we obtain that the partial $\Sigma^{\boldsymbol{\mathcal{A}}}$-algebra of path term classes $[\mathbf{PT}_{\boldsymbol{\mathcal{A}}}]$ is also a partial $\Sigma^{\boldsymbol{\mathcal{A}}}$-algebra in the variety $\mathcal{V}(\boldsymbol{\mathcal{E}}^{\boldsymbol{\mathcal{A}}})$ and isomorphic to $\mathbf{T}_{\boldsymbol{\mathcal{E}}^{\boldsymbol{\mathcal{A}}}}(\mathbf{Pth}_{\boldsymbol{\mathcal{A}}})$.

\section{
\texorpdfstring
{A $\mathrm{QE}$-variety of many-sorted partial $\Sigma^{\boldsymbol{\mathcal{A}}}$-algebras}
{A variety of partial algebras}
}

In this section we introduce, for the rewriting system $\boldsymbol{\mathcal{A}} = ((\mathbf{\Sigma},X),\mathcal{A})$, the $\mathrm{QE}$-variety $\mathcal{V}(\boldsymbol{\mathcal{E}}^{\boldsymbol{\mathcal{A}}})$ of the partial algebras relative to the categorial signature $\Sigma^{\boldsymbol{\mathcal{A}}}$ and we  investigate their connections to the partial $\Sigma^{\boldsymbol{\mathcal{A}}}$-algebras of path classes $[\mathbf{Pth}_{\boldsymbol{\mathcal{A}}}]$, introduced in Proposition~\ref{PCHCatAlg}, and by extension, in virtue of Theorem~\ref{TIso}, to the partial $\Sigma^{\boldsymbol{\mathcal{A}}}$-algebras of path term classes $[\mathbf{PT}_{\boldsymbol{\mathcal{A}}}]$, introduced in Proposition~\ref{PPTQCatAlg}.

\begin{restatable}{definition}{DVar}
\label{DVar}
\index{variety!first-order!$\boldsymbol{\mathcal{E}}^{\boldsymbol{\mathcal{A}}}$}
\index{variety!first-order!$\mathbf{Palg}(\boldsymbol{\mathcal{E}}^{\boldsymbol{\mathcal{A}}})$}
For the rewriting system $\boldsymbol{\mathcal{A}} = ((\mathbf{\Sigma},X),\mathcal{A})$,
we will denote by $((\mathbf{\Sigma}^{\boldsymbol{\mathcal{A}}},V^{S}),\mathcal{E}^{\boldsymbol{\mathcal{A}}})$, written $\boldsymbol{\mathcal{E}}^{\boldsymbol{\mathcal{A}}}$ for short, the specification in which 
$\mathbf{\Sigma}^{\boldsymbol{\mathcal{A}}}$ is $(S,\Sigma^{\boldsymbol{\mathcal{A}}})$, $V^{S}$ a fixed $S$-sorted set with a countable infinity of variables in each coordinate, and $\mathcal{E}^{\boldsymbol{\mathcal{A}}}$ the subset of $\mathrm{QE}(\Sigma^{\boldsymbol{\mathcal{A}}})_{V^{S}}$ (which, by Definition~\ref{DQEEq}, is the set of all $\mathrm{QE}$-equations with variables in $V^{S}$), consisting of the following equations:

For every $(\mathbf{s}, s)\in S^{\star}\times S$, every operation symbol $\sigma\in \Sigma_{\mathbf{s}, s}$, and every family of variables $(x_{j})_{j\in \bb{\mathbf{s}}}\in V^{S}_{\mathbf{s}}$, the operation $\sigma$ applied to the family $(x_{j})_{j\in \bb{\mathbf{s}}}$ is always defined. Formally, we have the following equation:
\allowdisplaybreaks
\begin{align*}\label{DVarA0}\tag{A0}
\sigma((x_{j})_{j\in \bb{\mathbf{s}}})\overset{\mathrm{e}}{=}\sigma((x_{j})_{j\in \bb{\mathbf{s}}}).
\end{align*}

For every sort $s\in S$ and every variable $x\in V^{S}_{s}$, the $0$-source and $0$-target of $x$ is always defined. Formally, we have the following equations:
\allowdisplaybreaks
\begin{align*}\label{DVarA1}\tag{A1}
\mathrm{sc}^{0}_{s}(x)&\overset{\mathrm{e}}{=}\mathrm{sc}^{0}_{s}(x);
&
\mathrm{tg}^{0}_{s}(x)&\overset{\mathrm{e}}{=}\mathrm{tg}^{0}_{s}(x).
\end{align*}

For every sort $s\in S$ and every variable $x\in V^{S}_{s}$, we have the following equations:
\allowdisplaybreaks
\begin{align*}\label{DVarA2}
\mathrm{sc}^{0}_{s}(\mathrm{sc}^{0}_{s}(x))&\overset{\mathrm{e}}{=}\mathrm{sc}^{0}_{s}(x);
&
\mathrm{sc}^{0}_{s}(\mathrm{tg}^{0}_{s}(x))&\overset{\mathrm{e}}{=}\mathrm{tg}^{0}_{s}(x);\\
\mathrm{tg}^{0}_{s}(\mathrm{sc}^{0}_{s}(x))&\overset{\mathrm{e}}{=}\mathrm{sc}^{0}_{s}(x);
&
\mathrm{tg}^{0}_{s}(\mathrm{tg}^{0}_{s}(x))&\overset{\mathrm{e}}{=}\mathrm{tg}^{0}_{s}(x).\tag{A2}
\end{align*}
In other words, $\mathrm{sc}^{0}_{s}$ and $\mathrm{tg}^{0}_{s}$ are right zeros. In particular, $\mathrm{sc}^{0}_{s}$ and $\mathrm{tg}^{0}_{s}$ are idempotent.

For every sort $s\in S$ and every pair of variables $x,y\in V^{S}_{s}$,  $x\circ^{0}_{s} y$ is defined if and only if the $0$-target of $y$ is equal to the $0$-source of $x$. Formally, we have the following conditional equations:
\allowdisplaybreaks
\begin{align*}\label{DVarA3}
x\circ^{0}_{s}y\overset{\mathrm{e}}{=}x\circ^{0}_{s}y \;\;
&\to \;\;
\mathrm{sc}^{0}_{s}(x)\overset{\mathrm{e}}{=}\mathrm{tg}^{0}_{s}(y);\\
\mathrm{sc}^{0}_{s}(x)\overset{\mathrm{e}}{=}\mathrm{tg}^{0}_{s}(y) \;\;
&\to \;\;
x\circ^{0}_{s}y\overset{\mathrm{e}}{=}x\circ^{0}_{s}y.\tag{A3}
\end{align*}

For every sort $s\in S$ and every pair of variables $x,y\in V^{S}_{s}$,  if $x\circ^{0}_{s} y$ is defined, then the $0$-source of $x\circ^{0}_{s}y$ is that of $y$ and the $0$-target of $x\circ^{0}_{s}y$ is that of $x$. Formally, we have the following conditional equations:
\allowdisplaybreaks
\begin{align*}\label{DVarA4}
x\circ^{0}_{s}y\overset{\mathrm{e}}{=}x\circ^{0}_{s}y \;\;
&\to \;\;
\mathrm{sc}^{0}_{s}(x\circ^{0}_{s}y)\overset{\mathrm{e}}{=}\mathrm{sc}^{0}_{s}(y);\\
x\circ^{0}_{s}y\overset{\mathrm{e}}{=}x\circ^{0}_{s}y \;\;
&\to \;\;
\mathrm{tg}^{0}_{s}(x\circ^{0}_{s}y)\overset{\mathrm{e}}{=}\mathrm{tg}^{0}_{s}(x).\tag{A4}
\end{align*}

For every sort $s\in S$ and every variable $x\in V^{S}_{s}$, the compositions $x\circ^{0}_{s}\mathrm{sc}^{0}_{s}(x)$ and $\mathrm{tg}^{0}_{s}(x)\circ^{0}_{s}x$ are always defined and are equal to $x$, i.e., $\mathrm{sc}^{0}_{s}(x)$ is a right unit element for the $0$-composition with $x$ and $\mathrm{tg}^{0}_{s}(x)$ is a left unit element for the $0$-composition with $x$. Formally, the following equations are satisfied
\allowdisplaybreaks
\begin{align*}\label{DVarA5}\tag{A5}
x\circ^{0}_{s}\mathrm{sc}^{0}_{s}(x)&\overset{\mathrm{e}}{=} x;
&
\mathrm{tg}^{0}_{s}(x)\circ^{0}_{s}x&\overset{\mathrm{e}}{=} x.
\end{align*}

For every sort $s\in S$ and every triple of variables $x, y, z\in V^{S}_{s}$, if the $0$-compositions $x\circ^{0}_{s}y$ and $y\circ^{0}_{s}z$ are defined, then the $0$-compositions $x\circ^{0}_{s}(y\circ^{0}_{s}z)$ and $(x\circ^{0}_{s}y)\circ^{0}_{s}z$ are defined and they are equal, i.e., the $0$-composition, when defined, is associative. Formally, we have the following conditional equation:
\allowdisplaybreaks
\begin{align*}\label{DVarA6}\tag{A6}
(x\circ^{0}_{s}y\overset{\mathrm{e}}{=} x\circ^{0}_{s}y)
\wedge
(y\circ^{0}_{s}z\overset{\mathrm{e}}{=} y\circ^{0}_{s}z)\;\;
&\to \;\;
(x\circ^{0}_{s}y)\circ^{0}_{s}z\overset{\mathrm{e}}{=} x\circ^{0}_{s}(y\circ^{0}_{s}z).
\end{align*}

For every $(\mathbf{s}, s)\in S^{\star}\times S$, every operation symbol $\sigma\in \Sigma_{\mathbf{s}, s}$, and every family of variables $(x_{j})_{j\in \bb{\mathbf{s}}}\in V^{S}_{\mathbf{s}}$, the $0$-source of $\sigma((x_{j})_{j\in \bb{\mathbf{s}}})$ is equal to $\sigma$ applied to the family $((\mathrm{sc}^{0}_{s_{j}}(x_{j}))_{j\in\bb{\mathbf{s}}})$, and the $0$-target of $\sigma((x_{j})_{j\in \bb{\mathbf{s}}})$ is equal to $\sigma$ applied to the family $((\mathrm{tg}^{0}_{s_{j}}(x_{j}))_{j\in\bb{\mathbf{s}}})$. Formally, we have the following equations:
\allowdisplaybreaks
\begin{align*}\label{DVarA7}
\mathrm{sc}^{0}_{s}(\sigma((x_{j})_{j\in \bb{\mathbf{s}}}))\overset{\mathrm{e}}{=}\sigma((\mathrm{sc}^{0}_{s_{j}}(x_{j}))_{j\in\bb{\mathbf{s}}});\\
\mathrm{tg}^{0}_{s}(\sigma((x_{j})_{j\in \bb{\mathbf{s}}}))\overset{\mathrm{e}}{=}\sigma((\mathrm{tg}^{0}_{s_{j}}(x_{j}))_{j\in\bb{\mathbf{s}}}).\tag{A7}
\end{align*}

For every $(\mathbf{s}, s)\in S^{\star}\times S$, every operation symbol $\sigma\in \Sigma_{\mathbf{s}, s}$, and every pair of families of variables $(x_{j})_{j\in \bb{\mathbf{s}}}, (y_{j})_{j\in \bb{\mathbf{s}}}\in V^{S}_{\mathbf{s}}$,  if, for every $j\in\bb{\mathbf{s}}$, the $0$-compositions $x_{j}\circ^{0}_{s_{j}}y_{j}$ are defined, then the $0$-composition $\sigma((x_{j})_{j\in \bb{\mathbf{s}}})\circ^{0}_{s}\sigma((y_{j})_{j\in \bb{\mathbf{s}}})$ is defined and it is equal to $\sigma$ applied to the family $(x_{j}\circ^{0}_{s_{j}}y_{j})_{j\in \bb{\mathbf{s}}}$. Formally, we have the following conditional equation:
\allowdisplaybreaks
\begin{multline*}\label{DVarA8}
\textstyle
\bigwedge_{j\in\bb{\mathbf{s}}}
(x_{j}\circ^{0}_{s_{j}}y_{j}\overset{\mathrm{e}}{=}x_{j}\circ^{0}_{s_{j}}y_{j})
\;\;\to\;\;\\
\sigma((x_{j}\circ^{0}_{s_{j}}y_{j})_{j\in \bb{\mathbf{s}}})
\overset{\mathrm{e}}{=}
\sigma((x_{j})_{j\in \bb{\mathbf{s}}})\circ^{0}_{s}\sigma((y_{j})_{j\in \bb{\mathbf{s}}})
\tag{A8}
\end{multline*}

For every sort $s\in S$ and every rewrite rule $\mathfrak{p}\in\mathcal{A}_{s}$, $\mathfrak{p}$ is always  defined. Formally, we have the following equation:
\allowdisplaybreaks
\begin{align*}\label{DVarA9}\tag{A9}
\mathfrak{p}\overset{\mathrm{e}}{=}\mathfrak{p}.
\end{align*}

We will let $\mathbf{PAlg}(\boldsymbol{\mathcal{E}}^{\boldsymbol{\mathcal{A}}})$ stand for the category canonically associated to the $\mathrm{QE}$-variety $\mathcal{V}(\boldsymbol{\mathcal{E}}^{\boldsymbol{\mathcal{A}}})$ determined by the specification $\boldsymbol{\mathcal{E}}^{\boldsymbol{\mathcal{A}}}$.
\end{restatable}

\begin{remark}
A model of axioms $A1$--$A6$ is an $S$-sorted category, see Definition~\ref{DnCat}. On the other hand, if one forgets the sorts involved in axioms $A1$--$A6$, which are $\mathrm{QE}$-equations, then a model of them is, simply, a category. We point out that Burmeister~\cite{pb02} states that the notion of category is definable by means of $\mathrm{ECE}$-equations and not only by means of $\mathrm{QE}$-equations.
\end{remark}

We next prove that the partial $\Sigma^{\boldsymbol{\mathcal{A}}}$-algebra $[\mathbf{Pth}_{\boldsymbol{\mathcal{A}}}]$ is a model of $\mathcal{E}^{\boldsymbol{\mathcal{A}}}$.

\begin{restatable}{proposition}{PPthVar}
\label{PPthVar}
$[\mathbf{Pth}_{\boldsymbol{\mathcal{A}}}]$ is a partial $\Sigma^{\boldsymbol{\mathcal{A}}}$-algebra in 
$\mathbf{PAlg}(\boldsymbol{\mathcal{E}}^{\boldsymbol{\mathcal{A}}})$.
\end{restatable}

\begin{proof}
We prove that all axioms defining $\mathbf{PAlg}(\boldsymbol{\mathcal{E}}^{\boldsymbol{\mathcal{A}}})$ are valid in $[\mathbf{Pth}_{\boldsymbol{\mathcal{A}}}]$. In this regard we recall from Proposition~\ref{PCHCatAlg} that, since $[\mathbf{Pth}_{\boldsymbol{\mathcal{A}}}]$ is a quotient of $\mathbf{Pth}_{\boldsymbol{\mathcal{A}}}$, the interpretation of the operation symbols in $[\mathbf{Pth}_{\boldsymbol{\mathcal{A}}}]$ is defined in terms of the interpretation of the operation symbols in $\mathbf{Pth}_{\boldsymbol{\mathcal{A}}}$. 

Axiom~\ref{DVarA0}. Let $(\mathbf{s},s)$ be an element of $S^{\star}\times S$, $\sigma$ an operation symbol in $\Sigma_{\mathbf{s},s}$ and $([\mathfrak{P}_{j}]^{}_{s_{j}})_{j\in\bb{\mathbf{s}}}$ a family of path classes in $[\mathrm{Pth}_{\boldsymbol{\mathcal{A}}}]_{\mathbf{s}}$. Then, by Proposition~\ref{PCHCatAlg}, we have that the interpretation of $\sigma$ as an operation symbol in the many-sorted partial $\Sigma^{\boldsymbol{\mathcal{A}}}$-algebra $[\mathbf{Pth}_{\boldsymbol{\mathcal{A}}}]$ is always defined.

This proves that Axiom~\ref{DVarA0} is valid in $[\mathbf{Pth}_{\boldsymbol{\mathcal{A}}}]$.

Axiom~\ref{DVarA1}. We recall from Proposition~\ref{PPthCatAlg} that, for every sort $s\in S$, the $0$-source and $0$-target operations are totally defined in the many-sorted partial $\Sigma^{\boldsymbol{\mathcal{A}}}$-algebra $\mathbf{Pth}_{\boldsymbol{\mathcal{A}}}$. Let  $s$ be a sort in $S$ and $[\mathfrak{P}]^{}_{s}$ a path class in $[\mathrm{Pth}_{\boldsymbol{\mathcal{A}}}]_{s}$, then
\allowdisplaybreaks
\begin{align*}
\mathrm{sc}_{s}
^{0[\mathbf{Pth}_{\boldsymbol{\mathcal{A}}}]}
\left(
\left[\mathfrak{P}
\right]^{}_{s}
\right)
&=
\left[
\mathrm{sc}_{s}
^{0\mathbf{Pth}_{\boldsymbol{\mathcal{A}}}}
\left(
\mathfrak{P}
\right)
\right]^{}_{s};
\\
\mathrm{tg}_{s}
^{0[\mathbf{Pth}_{\boldsymbol{\mathcal{A}}}]}
\left(
\left[\mathfrak{P}
\right]^{}_{s}
\right)
&=
\left[
\mathrm{tg}_{s}
^{0\mathbf{Pth}_{\boldsymbol{\mathcal{A}}}}
\left(
\mathfrak{P}
\right)
\right]^{}_{s}.
\end{align*}

This proves that Axiom~\ref{DVarA1} is valid in $[\mathbf{Pth}_{\boldsymbol{\mathcal{A}}}]$.

Axiom~\ref{DVarA2}. Let $s$ be a sort in $S$ and $[\mathfrak{P}]^{}_{s}$ a path class in $[\mathrm{Pth}_{\boldsymbol{\mathcal{A}}}]_{s}$. Then, by  Proposition~\ref{PCHVarA2}, we have 
\allowdisplaybreaks
\begin{align*}
\mathrm{sc}^{0[\mathbf{Pth}_{\boldsymbol{\mathcal{A}}}]}_{s}\left(
\mathrm{sc}^{0[\mathbf{Pth}_{\boldsymbol{\mathcal{A}}}]}_{s}\left(
[\mathfrak{P}]_{s}
\right)\right)
&=
\mathrm{sc}^{0[\mathbf{Pth}_{\boldsymbol{\mathcal{A}}}]}_{s}\left(
[\mathfrak{P}]_{s}
\right);
\\
\mathrm{sc}^{0[\mathbf{Pth}_{\boldsymbol{\mathcal{A}}}]}_{s}\left(
\mathrm{tg}^{0[\mathbf{Pth}_{\boldsymbol{\mathcal{A}}}]}_{s}\left(
[\mathfrak{P}]_{s}
\right)\right)
&=
\mathrm{tg}^{0[\mathbf{Pth}_{\boldsymbol{\mathcal{A}}}]}_{s}\left(
[\mathfrak{P}]_{s}
\right);
\\
\mathrm{tg}^{0[\mathbf{Pth}_{\boldsymbol{\mathcal{A}}}]}_{s}\left(
\mathrm{sc}^{0[\mathbf{Pth}_{\boldsymbol{\mathcal{A}}}]}_{s}\left(
[\mathfrak{P}]_{s}
\right)\right)
&=
\mathrm{sc}^{0[\mathbf{Pth}_{\boldsymbol{\mathcal{A}}}]}_{s}\left(
[\mathfrak{P}]_{s}
\right);
\\
\mathrm{tg}^{0[\mathbf{Pth}_{\boldsymbol{\mathcal{A}}}]}_{s}\left(
\mathrm{tg}^{0[\mathbf{Pth}_{\boldsymbol{\mathcal{A}}}]}_{s}\left(
[\mathfrak{P}]_{s}
\right)\right)
&=
\mathrm{tg}^{0[\mathbf{Pth}_{\boldsymbol{\mathcal{A}}}]}_{s}\left(
[\mathfrak{P}]_{s}
\right).
\end{align*}

This proves that Axiom~\ref{DVarA2} is valid in $[\mathbf{Pth}_{\boldsymbol{\mathcal{A}}}]$.

Axiom~\ref{DVarA3}. Let $s$ be a sort in $S$ and $\mathfrak{Q}, \mathfrak{P}$ a pair of paths in $\mathrm{Pth}_{\boldsymbol{\mathcal{A}},s}$. By Proposition~\ref{PCHVarA3}, we can affirm that the following conditions are equivalent
\begin{itemize}
\item[(i)] $[\mathfrak{Q}]^{}_{s}
\circ_{s}^{0[\mathbf{Pth}_{\boldsymbol{\mathcal{A}}}]}
[\mathfrak{P}]^{}_{s}$ is defined;
\item[(ii)] $\mathrm{sc}_{s}
^{0[\mathbf{Pth}_{\boldsymbol{\mathcal{A}}}]}
(\mathfrak{[Q}]^{}_{s})
=
\mathrm{tg}_{s}
^{0[\mathbf{Pth}_{\boldsymbol{\mathcal{A}}}]}
(\mathfrak{[P}]^{}_{s})$.
\end{itemize}

This proves that Axiom~\ref{DVarA3} is valid in $[\mathbf{Pth}_{\boldsymbol{\mathcal{A}}}]$.

Axiom~\ref{DVarA4}. Let $s$ be a sort in $S$ and let $[\mathfrak{Q}]_{s}$ and $[\mathfrak{P}]^{}_{s}$ be a pair of path classes in $[\mathrm{Pth}_{\boldsymbol{\mathcal{A}}}]_{s}$ such that the $0$-composition 
$[\mathfrak{Q}]^{}_{s}
\circ_{s}^{0[\mathbf{Pth}_{\boldsymbol{\mathcal{A}}}]}
[\mathfrak{P}]^{}_{s}$
is defined.  Then, by Proposition~\ref{PCHVarA4}, we have 
\allowdisplaybreaks
\begin{align*}
\mathrm{sc}^{0[\mathbf{Pth}_{\boldsymbol{\mathcal{A}}}]}_{s}\left(
[\mathfrak{Q}]_{s}
\circ^{0[\mathbf{Pth}_{\boldsymbol{\mathcal{A}}}]}_{s}
[\mathfrak{P}]_{s}
\right)
&=
\mathrm{sc}^{0[\mathbf{Pth}_{\boldsymbol{\mathcal{A}}}]}_{s}\left(
[\mathfrak{P}]_{s}
\right);
\\
\mathrm{tg}^{0[\mathbf{Pth}_{\boldsymbol{\mathcal{A}}}]}_{s}\left(
[\mathfrak{Q}]_{s}
\circ^{0[\mathbf{Pth}_{\boldsymbol{\mathcal{A}}}]}_{s}
[\mathfrak{P}]_{s}
\right)
&=
\mathrm{tg}^{0[\mathbf{Pth}_{\boldsymbol{\mathcal{A}}}]}_{s}\left(
[\mathfrak{Q}]_{s}
\right).
\end{align*}

This proves that Axiom~\ref{DVarA4} is valid in $[\mathbf{Pth}_{\boldsymbol{\mathcal{A}}}]$.

Axiom~\ref{DVarA5}. Let $s$ be a sort in $S$ and $[\mathfrak{P}]_{s}$ an equivalence class in $[\mathrm{Pth}_{\boldsymbol{\mathcal{A}}}]_{s}$. Then, according to Proposition~\ref{PCHVarA5}, we have that 
\allowdisplaybreaks
\begin{align*}
[\mathfrak{P}]_{s}
\circ^{0[\mathbf{Pth}_{\boldsymbol{\mathcal{A}}}]}_{s}
\left(
\mathrm{sc}^{0[\mathbf{Pth}_{\boldsymbol{\mathcal{A}}}]}_{s}\left(
[\mathfrak{P}]_{s}
\right)
\right)
&=
[\mathfrak{P}]_{s};
\\
\left(\mathrm{tg}^{0[\mathbf{Pth}_{\boldsymbol{\mathcal{A}}}]}_{s}\left(
[\mathfrak{P}]_{s}
\right)
\right)
\circ^{0[\mathbf{Pth}_{\boldsymbol{\mathcal{A}}}]}_{s}
[\mathfrak{P}]_{s}
&=
[\mathfrak{P}]_{s}.
\end{align*}

This proves that Axiom~\ref{DVarA5} is valid in $[\mathbf{Pth}_{\boldsymbol{\mathcal{A}}}]$.

Axiom~\ref{DVarA6}. Let $s$ be a sort in $S$ and let $[\mathfrak{P}]^{}_{s}$, $[\mathfrak{Q}]^{}_{s}$, and $[\mathfrak{R}]_{s}^{}$ be path classes in $[\mathrm{Pth}_{\boldsymbol{\mathcal{A}}}]_{s}$ such that 
$
[\mathfrak{Q}]^{}_{s}
\circ^{0[\mathbf{Pth}_{\boldsymbol{\mathcal{A}}}]}_{s}
[\mathfrak{P}]^{}_{s}
$ 
and 
$
[\mathfrak{R}]^{}_{s}
\circ^{0[\mathbf{Pth}_{\boldsymbol{\mathcal{A}}}]}_{s}
[\mathfrak{Q}]^{}_{s}
$
are defined. Then, according to Proposition~\ref{PCHVarA6}, we have that 
\[
[\mathfrak{R}]_{s}
\circ^{0[\mathbf{Pth}_{\boldsymbol{\mathcal{A}}}]}_{s}\left(
[\mathfrak{Q}]_{s}
\circ^{0[\mathbf{Pth}_{\boldsymbol{\mathcal{A}}}]}_{s}
[\mathfrak{P}]_{s}
\right)
=
\left(
[\mathfrak{R}]_{s}
\circ^{0[\mathbf{Pth}_{\boldsymbol{\mathcal{A}}}]}_{s}
[\mathfrak{Q}]_{s}
\right)
\circ^{0[\mathbf{Pth}_{\boldsymbol{\mathcal{A}}}]}_{s}
[\mathfrak{P}]_{s}
.
\]

This proves that Axiom~\ref{DVarA6} is valid in $[\mathbf{Pth}_{\boldsymbol{\mathcal{A}}}]$.

Axiom~\ref{DVarA7}. Let $(\mathbf{s},s)$ be an element of $S^{\star}\times S$, $\sigma$ an operation symbol in $\Sigma_{\mathbf{s},s}$ and $([\mathfrak{P}_{j}]^{}_{s_{j}})_{j\in\bb{\mathbf{s}}}$ a family of path classes in $[\mathrm{Pth}_{\boldsymbol{\mathcal{A}}}]_{\mathbf{s}}$. 
Then, according to Proposition~\ref{PCHVarA7}, we have that 
\allowdisplaybreaks
\begin{align*}
\mathrm{sc}^{0[\mathbf{Pth}_{\boldsymbol{\mathcal{A}}}]}_{s}\left(
\sigma^{[\mathbf{Pth}_{\boldsymbol{\mathcal{A}}}]}
\left(\left(
\left[
\mathfrak{P}_{j}
\right]_{s_{j}}
\right)_{j\in\bb{\mathbf{s}}}
\right)\right)
&=
\sigma^{[\mathbf{Pth}_{\boldsymbol{\mathcal{A}}}]}
\left(\left(
\mathrm{sc}^{0[\mathbf{Pth}_{\boldsymbol{\mathcal{A}}}]}_{s_{j}}\left(
\left[
\mathfrak{P}_{j}
\right]_{s_{j}}
\right)\right)_{j\in\bb{\mathbf{s}}}\right);
\\
\mathrm{tg}^{0[\mathbf{Pth}_{\boldsymbol{\mathcal{A}}}]}_{s}\left(
\sigma^{[\mathbf{Pth}_{\boldsymbol{\mathcal{A}}}]}
\left(\left(
\left[
\mathfrak{P}_{j}
\right]_{s_{j}}
\right)_{j\in\bb{\mathbf{s}}}
\right)\right)
&=
\sigma^{[\mathbf{Pth}_{\boldsymbol{\mathcal{A}}}]}
\left(\left(
\mathrm{tg}^{0[\mathbf{Pth}_{\boldsymbol{\mathcal{A}}}]}_{s_{j}}\left(
\left[
\mathfrak{P}_{j}
\right]_{s_{j}}
\right)\right)_{j\in\bb{\mathbf{s}}}\right).
\end{align*}

This proves that Axiom~\ref{DVarA7} is valid in $[\mathbf{Pth}_{\boldsymbol{\mathcal{A}}}]$.

Axiom~\ref{DVarA8}.  Let $(\mathbf{s},s)$ be an element of $S^{\star}\times S$, $\sigma\in\Sigma_{\mathbf{s},s}$ and let $([\mathfrak{P}_{j}]^{}_{s_{j}})_{j\in\bb{\mathbf{s}}}$ and $([\mathfrak{Q}_{j}]^{}_{s_{j}})_{j\in\bb{\mathbf{s}}}$ be two families of path classes in $[\mathrm{Pth}_{\boldsymbol{\mathcal{A}}}]_{\mathbf{s}}$ such that, for every $j\in\bb{\mathbf{s}}$, the $0$-compositions
$
[\mathfrak{Q}_{j}]^{}_{s_{j}}
\circ_{s_{j}}^{0[\mathbf{Pth}_{\boldsymbol{\mathcal{A}}}]}
[\mathfrak{P}_{j}]^{}_{s_{j}}
$
are defined. 

Then, according to Proposition~\ref{PCHVarA8}, we have that
\allowdisplaybreaks
\begin{multline*}
\sigma^{[\mathbf{Pth}_{\boldsymbol{\mathcal{A}}}]}
\left(\left(
\left[\mathfrak{Q}_{j}\right]_{s_{j}}
\circ^{0[\mathbf{Pth}_{\boldsymbol{\mathcal{A}}}]}_{s_{j}}
\left[\mathfrak{P}_{j}\right]_{s_{j}}
\right)_{j\in\bb{\mathbf{s}}}
\right)
\\=
\sigma^{[\mathbf{Pth}_{\boldsymbol{\mathcal{A}}}]}
\left(\left(
\left[\mathfrak{Q}_{j}\right]_{s_{j}}
\right)_{j\in\bb{\mathbf{s}}}
\right)
\circ^{0[\mathbf{Pth}_{\boldsymbol{\mathcal{A}}}]}_{s}
\sigma^{[\mathbf{Pth}_{\boldsymbol{\mathcal{A}}}]}
\left(\left(
\left[\mathfrak{P}_{j}\right]_{s_{j}}
\right)_{j\in\bb{\mathbf{s}}}
\right).
\end{multline*}

This proves that Axiom~\ref{DVarA8} is valid in $[\mathbf{Pth}_{\boldsymbol{\mathcal{A}}}]$.

Axiom~\ref{DVarA9}. Let $s$ be a sort in $S$ and $\mathfrak{p}$ a rewrite rule in $\mathcal{A}_{s}$. Then, according to Proposition~\ref{PCHCatAlg} and Definition~\ref{DCHEch}, we have that the interpretation of $\mathfrak{p}$ as a constant in $[\mathbf{Pth}_{\boldsymbol{\mathcal{A}}}]$ is given by $\mathrm{ech}^{([1],\mathcal{A})}_{s}(\mathfrak{p})$, that is, the equivalence  class under the kernel of the Curry-Howard mapping of the echelon determined by $\mathfrak{p}$.

This proves that Axiom~\ref{DVarA9} is valid in $[\mathbf{Pth}_{\boldsymbol{\mathcal{A}}}]$.

This completes the proof.
\end{proof}

\begin{remark} Among the axioms in $\mathcal{E}^{\boldsymbol{\mathcal{A}}}$, the Axiom~\ref{DVarA8} is the only one that is not valid in the many-sorted partial $\Sigma^{\boldsymbol{\mathcal{A}}}$-algebra $\mathbf{Pth}_{\boldsymbol{\mathcal{A}}}$ in case the signature $\Sigma$ contains non-constant operation symbols.
\end{remark}

\begin{restatable}{corollary}{CPTVar}
\label{CPTVar} $[\mathbf{PT}_{\boldsymbol{\mathcal{A}}}]$ is a many-sorted partial  $\Sigma^{\boldsymbol{\mathcal{A}}}$-algebra in $\mathbf{PAlg}(\boldsymbol{\mathcal{E}}^{\boldsymbol{\mathcal{A}}})$.
\end{restatable}
\begin{proof}
It follows from Proposition~\ref{PPthVar} and Theorem~\ref{TIso}.
\end{proof}

\section{
\texorpdfstring
{Freedom in $\mathcal{V}(\boldsymbol{\mathcal{E}}^{\boldsymbol{\mathcal{A}}})$}
{Freedom}
}

The aim of this section is to prove that the partial $\Sigma^{\boldsymbol{\mathcal{A}}}$-algebra of path classes $[\mathbf{Pth}_{\boldsymbol{\mathcal{A}}}]$ is a $\mathbf{PAlg}(\boldsymbol{\mathcal{E}}^{\boldsymbol{\mathcal{A}}})$-universal solution of the many-sorted partial $\Sigma^{\boldsymbol{\mathcal{A}}}$-algebra $\mathbf{Pth}_{\boldsymbol{\mathcal{A}}}$. In virtue of Theorem~\ref{TIso}, the same condition will apply to the partial $\Sigma^{\boldsymbol{\mathcal{A}}}$-algebra of path term classes $[\mathbf{PT}_{\boldsymbol{\mathcal{A}}}]$.

The reader is advised to consult Theorem~\ref{TFreeAdj} for a description of the free many-sorted  partial $\Sigma^{\boldsymbol{\mathcal{A}}}$-algebra associated to a $\mathrm{QE}$-variety of many-sorted partial $\Sigma^{\boldsymbol{\mathcal{A}}}$-algebras relative to a many-sorted partial $\Sigma^{\boldsymbol{\mathcal{A}}}$-algebra. 

\begin{restatable}{definition}{DVarAbv}
\label{DVarAbv}
\index{variety!first-order!$\mathbf{T}_{\boldsymbol{\mathcal{E}}^{\boldsymbol{\mathcal{A}}}}
(\mathbf{Pth}_{\boldsymbol{\mathcal{A}}})$}
Taking into account Theorem~\ref{TFreeAdj}, we introduce the following  abbreviations
\allowdisplaybreaks
\begin{align*}
\mathbf{Sch}
_{\boldsymbol{\mathcal{E}}^{\boldsymbol{\mathcal{A}}}}
(\mathbf{Pth}_{\boldsymbol{\mathcal{A}}})
&=
\textstyle
\bigcap
_{
f\in
\bigcup_{
\mathbf{B}\in
\mathbf{PAlg}
(
\boldsymbol{\mathcal{E}}^{\boldsymbol{\mathcal{A}}}
)
}
\mathrm{Hom}(
\mathbf{Pth}_{\boldsymbol{\mathcal{A}}}
,
\mathbf{B})
}
\mathbf{Sch}(f)
,
\\
\equiv^{[1]}
&=
\textstyle
\bigcap_{f\in\bigcup_{
\mathbf{B}\in
\mathbf{PAlg}
(
\boldsymbol{\mathcal{E}}^{\boldsymbol{\mathcal{A}}}
)
}\mathrm{Hom}(
\mathbf{Pth}_{\boldsymbol{\mathcal{A}}}
,
\mathbf{B}
)
}
\mathrm{SKer}(f)
,\,\text{and}
\\
\mathbf{T}_{\boldsymbol{\mathcal{E}}^{\boldsymbol{\mathcal{A}}}}
(\mathbf{Pth}_{\boldsymbol{\mathcal{A}}})
&=
\mathbf{Sch}
_{\boldsymbol{\mathcal{E}}^{\boldsymbol{\mathcal{A}}}}
(\mathbf{Pth}_{\boldsymbol{\mathcal{A}}})
/
{
\equiv^{[1]}
}
,
\end{align*}
where, we recall, $\mathbf{Sch}(f)$, the Schmidt closed $\mathbf{Pth}_{\boldsymbol{\mathcal{A}}}$-initial extension of $f$, and $\mathrm{SKer}(f)$, the Schmidt Kernel of $f$, where defined in Proposition \ref{PSch}. 
\end{restatable}

The aim of this section is to prove:
$$
[\mathbf{Pth}_{\boldsymbol{\mathcal{A}}}]
\cong
\mathbf{T}_{\boldsymbol{\mathcal{E}}^{\boldsymbol{\mathcal{A}}}}
\left(\mathbf{Pth}_{\boldsymbol{\mathcal{A}}}\right).
$$

This proof will require several intermediate steps. We will start by introducing the following terminology.

\begin{restatable}{definition}{DVarAp}
\label{DVarAp}
We will denote by 
\begin{enumerate}
\item $\mathrm{pr}^{\equiv^{[1]}}$ the projection from $\mathrm{Sch}_{\boldsymbol{\mathcal{E}}^{\boldsymbol{\mathcal{A}}}}(\mathbf{Pth}_{\boldsymbol{\mathcal{A}}})$ to its quotient $\mathrm{T}_{\boldsymbol{\mathcal{E}}^{\boldsymbol{\mathcal{A}}}}(\mathbf{Pth}_{\boldsymbol{\mathcal{A}}})$;
\item  $\eta^{([1],\mathbf{Pth}_{\boldsymbol{\mathcal{A}}})}$ the $S$-sorted mapping from $\mathrm{Pth}_{\boldsymbol{\mathcal{A}}}$ to $\mathrm{T}_{\boldsymbol{\mathcal{E}}^{\boldsymbol{\mathcal{A}}}}(\mathbf{Pth}_{\boldsymbol{\mathcal{A}}})$ given by the composition
$
\eta^{([1],\mathbf{Pth}_{\boldsymbol{\mathcal{A}}})}=
\mathrm{pr}^{\equiv^{[1]}}\circ\eta^{(1,\mathbf{Pth}_{\boldsymbol{\mathcal{A}}})},
$
where $\eta^{(1,\mathbf{Pth}_{\boldsymbol{\mathcal{A}}})}$ is the insertion of generators, from $\mathrm{Pth}_{\boldsymbol{\mathcal{A}}}$ to $\mathrm{F}_{\Sigma^{\boldsymbol{\mathcal{A}}}}(\mathbf{Pth}_{\boldsymbol{\mathcal{A}}})$, introduced in Definition~\ref{DFP}.
\end{enumerate}
\end{restatable}

\begin{remark}\label{RVar} Let us recall from Proposition~\ref{PFreeComp} and Definition~\ref{DFP} that $\eta^{(1,\mathbf{Pth}_{\boldsymbol{\mathcal{A}}})}$ is a $\Sigma^{\boldsymbol{\mathcal{A}}}$-homomorphism from $\mathbf{Pth}_{\boldsymbol{\mathcal{A}}}$ to $\mathbf{F}_{\Sigma^{\boldsymbol{\mathcal{A}}}}(\mathbf{Pth}_{\boldsymbol{\mathcal{A}}})$. Moreover, since every Schmidt algebra involved in the construction of $\mathbf{Sch}_{\boldsymbol{\mathcal{E}}^{\boldsymbol{\mathcal{A}}}}(\mathbf{Pth}_{\boldsymbol{\mathcal{A}}})$ is $\mathbf{Pth}_{\boldsymbol{\mathcal{A}}}$-generated, we have that $\eta^{(1,\mathbf{Pth}_{\boldsymbol{\mathcal{A}}})}$ corestricts to $\mathbf{Sch}_{\boldsymbol{\mathcal{E}}^{\boldsymbol{\mathcal{A}}}}(\mathbf{Pth}_{\boldsymbol{\mathcal{A}}})$.
$$
\eta^{(1,\mathbf{Pth}_{\boldsymbol{\mathcal{A}}})}
\colon
\mathbf{Pth}_{\boldsymbol{\mathcal{A}}}
\mor
\mathbf{Sch}_{\boldsymbol{\mathcal{E}}^{\boldsymbol{\mathcal{A}}}}(\mathbf{Pth}_{\boldsymbol{\mathcal{A}}}).
$$

Moreover, taking into account Definition~\ref{DVarAbv}, the projection $\mathrm{pr}^{\equiv^{[1]}}$ is a $\Sigma^{\boldsymbol{\mathcal{A}}}$-homomorphism from $\mathbf{Sch}_{\boldsymbol{\mathcal{E}}^{\boldsymbol{\mathcal{A}}}}(\mathbf{Pth}_{\boldsymbol{\mathcal{A}}})$ to $\mathbf{T}_{\boldsymbol{\mathcal{E}}^{\boldsymbol{\mathcal{A}}}}(\mathbf{Pth}_{\boldsymbol{\mathcal{A}}})$
$$
\mathrm{pr}^{\equiv^{[1]}}
\colon
\mathbf{Sch}_{\boldsymbol{\mathcal{E}}^{\boldsymbol{\mathcal{A}}}}(\mathbf{Pth}_{\boldsymbol{\mathcal{A}}})
\mor
\mathbf{T}_{\boldsymbol{\mathcal{E}}^{\boldsymbol{\mathcal{A}}}}(\mathbf{Pth}_{\boldsymbol{\mathcal{A}}}).
$$

Taking into account these considerations and Definition~\ref{DVarAp}, $\eta^{([1],\mathbf{Pth}_{\boldsymbol{\mathcal{A}}})}$ is a  $\Sigma^{\boldsymbol{\mathcal{A}}}$-homomorphism from $\mathbf{Pth}_{\boldsymbol{\mathcal{A}}}$ to $\mathbf{T}_{\boldsymbol{\mathcal{E}}^{\boldsymbol{\mathcal{A}}}}(\mathbf{Pth}_{\boldsymbol{\mathcal{A}}})$
$$
\eta^{([1],\mathbf{Pth}_{\boldsymbol{\mathcal{A}}})}
\colon
\mathbf{Pth}_{\boldsymbol{\mathcal{A}}}
\mor
\mathbf{T}_{\boldsymbol{\mathcal{E}}^{\boldsymbol{\mathcal{A}}}}(\mathbf{Pth}_{\boldsymbol{\mathcal{A}}}).
$$
\end{remark}

The following result is a consequence of Proposition~\ref{PPthVar}. It states that the canonical projection from $\mathbf{Pth}_{\boldsymbol{\mathcal{A}}}$ to $[\mathbf{Pth}_{\boldsymbol{\mathcal{A}}}]$ can be extended  to a $\Sigma^{\boldsymbol{\mathcal{A}}}$-epimorphism from $\mathbf{T}_{\boldsymbol{\mathcal{E}}^{\boldsymbol{\mathcal{A}}}}
(\mathbf{Pth}_{\boldsymbol{\mathcal{A}}})$ to $[\mathbf{Pth}_{\boldsymbol{\mathcal{A}}}]$.

\begin{restatable}{corollary}{CVarPr}
\label{CVarPr}
For the $\Sigma^{\boldsymbol{\mathcal{A}}}$-epimorphism
$
\textstyle
\mathrm{pr}^{\mathrm{Ker}(\mathrm{CH}^{(1)})}
\colon
\mathbf{Pth}_{\boldsymbol{\mathcal{A}}}
\mor
[\mathbf{Pth}_{\boldsymbol{\mathcal{A}}}]
$, there exists a unique $\Sigma^{\boldsymbol{\mathcal{A}}}$-epimorphism
$$
\textstyle
\mathrm{pr}^{\mathrm{Ker}(\mathrm{CH}^{(1)})\mathsf{p}}
\colon
\mathbf{T}_{\boldsymbol{\mathcal{E}}^{\boldsymbol{\mathcal{A}}}}
(\mathbf{Pth}_{\boldsymbol{\mathcal{A}}})
\mor
[\mathbf{Pth}_{\boldsymbol{\mathcal{A}}}]
$$
such that $\mathrm{pr}^{\mathrm{Ker}(\mathrm{CH}^{(1)})\mathsf{p}}\circ\eta^{(1,\mathbf{Pth}_{\boldsymbol{\mathcal{A}}})}=\mathrm{pr}^{\mathrm{Ker}(\mathrm{CH}^{(1)})}$, i.e., such that
the diagram in Figure~\ref{FVarPr} commutes.
\end{restatable}
\begin{figure}
\begin{tikzpicture}
[ACliment/.style={-{To [angle'=45, length=5.75pt, width=4pt, round]},
},scale=0.8]
\node[] (1) at (0,0) [] {$\mathbf{Pth}_{\boldsymbol{\mathcal{A}}}$};
\node[] (3) at (6,-3) []
{$[\mathbf{Pth}_{\boldsymbol{\mathcal{A}}}]$};
\node[] (4) at (6,0) []
{$\mathbf{T}_{\boldsymbol{\mathcal{E}}^{\boldsymbol{\mathcal{A}}}}
(\mathbf{Pth}_{\boldsymbol{\mathcal{A}}})$};
\draw[ACliment, bend right=10]  (1) to node [below left]
{$\mathrm{pr}^{\mathrm{Ker}(\mathrm{CH}^{(1)})}$} (3);
\draw[ACliment]  (4) to node [right]
{$\mathrm{pr}^{\mathrm{Ker}(\mathrm{CH}^{(1)})\mathsf{p}}$} (3);
\draw[ACliment]  (1) to node [above]
{$\eta^{([1],\mathbf{Pth}_{\boldsymbol{\mathcal{A}}})}$} (4);
\end{tikzpicture}
\caption{The universal extension of $\mathrm{pr}^{\mathrm{Ker(CH)}}$.}\label{FVarPr}
\end{figure}

\begin{proof}
Let us note that $\mathrm{pr}^{\mathrm{Ker}(\mathrm{CH}^{(1)})}$ is a $\Sigma^{\boldsymbol{\mathcal{A}}}$-epimorphism from
$\mathbf{Pth}_{\boldsymbol{\mathcal{A}}}$ to $[\mathbf{Pth}_{\boldsymbol{\mathcal{A}}}]$ which, by Proposition~\ref{PPthVar}, is a many-sorted partial $\Sigma^{\boldsymbol{\mathcal{A}}}$-algebra in the $\mathrm{QE}$-variety $\mathbf{PAlg}(\boldsymbol{\mathcal{E}}^{\boldsymbol{\mathcal{A}}})$. Then, by the universal property of $\mathbf{T}_{\boldsymbol{\mathcal{E}}^{\boldsymbol{\mathcal{A}}}}
(\mathbf{Pth}_{\boldsymbol{\mathcal{A}}})$, Theorem~\ref{TFreeAdj}, there exists a unique $\Sigma^{\boldsymbol{\mathcal{A}}}$-homomorphism
$$
\mathrm{pr}^{\mathrm{Ker}(\mathrm{CH}^{(1)})\mathsf{p}}
\colon
\mathbf{T}_{\boldsymbol{\mathcal{E}}^{\boldsymbol{\mathcal{A}}}}
(\mathbf{Pth}_{\boldsymbol{\mathcal{A}}})
\mor
[\mathbf{Pth}_{\boldsymbol{\mathcal{A}}}]
$$
such that $\mathrm{pr}^{\mathrm{Ker}(\mathrm{CH}^{(1)})\mathsf{p}}\circ\eta^{([1],\mathbf{Pth}_{\boldsymbol{\mathcal{A}}})}=\mathrm{pr}^{\mathrm{Ker}(\mathrm{CH}^{(1)})}$. Moreover, since $\mathrm{pr}^{\mathrm{Ker}(\mathrm{CH}^{(1)})}$ is a $\Sigma^{\boldsymbol{\mathcal{A}}}$-epimorphism, we have that 
$\mathrm{pr}^{\mathrm{Ker}(\mathrm{CH}^{(1)})\mathsf{p}}$ is also a $\Sigma^{\boldsymbol{\mathcal{A}}}$-epimorphism.
\end{proof}

We now work towards the proof of the existence of an inverse of the $\Sigma^{\boldsymbol{\mathcal{A}}}$-homo\-mor\-phism $\mathrm{pr}^{\mathrm{Ker}(\mathrm{CH}^{(1)})\mathsf{p}}$. 

\begin{restatable}{proposition}{PVarKer}
\label{PVarKer}
The mapping 
$$
\mathrm{pr}^{
\equiv^{[1]}
}\circ\mathrm{ip}^{(1,X)@}\circ\mathrm{CH}^{(1)}
\colon
\mathbf{Pth}_{\boldsymbol{\mathcal{A}}}
\mor
\mathbf{T}_{\boldsymbol{\mathcal{E}}^{\boldsymbol{\mathcal{A}}}}
(
\mathbf{Pth}_{\boldsymbol{\mathcal{A}}}
)
$$
is a $\Sigma^{\boldsymbol{\mathcal{A}}}$-homomorphism such that
\[
\mathrm{Ker}(\mathrm{CH}^{(1)})
\subseteq \mathrm{Ker}(\mathrm{pr}^{
\equiv^{[1]}
}\circ\mathrm{ip}^{(1,X)@}\circ\mathrm{CH}^{(1)}).\]
\end{restatable}

\begin{proof}
Let us note that since every Schmidt algebra involved in the construction of 
$\mathbf{Sch}_{\boldsymbol{\mathcal{E}}^{\boldsymbol{\mathcal{A}}}}(\mathbf{Pth}_{\boldsymbol{\mathcal{A}}})$
is $\mathbf{Pth}_{\boldsymbol{\mathcal{A}}}$-generated and taking into account Proposition~\ref{PIpCH} we can corestrict the composition $\mathrm{ip}^{(1,X)@}\circ \mathrm{CH}^{(1)}$ to 
$\mathrm{Sch}_{\boldsymbol{\mathcal{E}}^{\boldsymbol{\mathcal{A}}}}(\mathbf{Pth}_{\boldsymbol{\mathcal{A}}})$, that is, 
\[
\mathrm{ip}^{(1,X)@}\circ \mathrm{CH}^{(1)}
\colon 
\mathrm{Pth}_{\boldsymbol{\mathcal{A}}}
\mor
\mathrm{Sch}_{\boldsymbol{\mathcal{E}}^{\boldsymbol{\mathcal{A}}}}(\mathbf{Pth}_{\boldsymbol{\mathcal{A}}}).
\]

This correstriction justifies  that  $\mathrm{pr}^{
\equiv^{[1]}
}\circ\mathrm{ip}^{(1,X)@}\circ\mathrm{CH}^{(1)}$ is well-defined.

We prove that $\mathrm{pr}^{
\equiv^{[1]}
}\circ\mathrm{ip}^{(1,X)@}\circ\mathrm{CH}^{(1)}$ is a $\Sigma^{\boldsymbol{\mathcal{A}}}$-homomorphism by checking the compatibility with all the different operations in $\Sigma^{\boldsymbol{\mathcal{A}}}$.

\textsf{Compatibility with the operations from $\Sigma$.}

Let $(\mathbf{s},s)$ be a pair in $\mathbf{S}^{\star}\times S$ and let $\sigma$ be an operation symbol in $\Sigma_{\mathbf{s},s}$, let $(\mathfrak{P}_{j})_{j\in\bb{\mathbf{s}}}$ be a family of paths in $\mathrm{Pth}_{\boldsymbol{\mathcal{A}},\mathbf{s}}$. 

The following chain of equalities holds
\begin{flushleft}
$\mathrm{pr}^{\equiv^{[1]}}_{s}\left(
\mathrm{ip}^{(1,X)@}_{s}\left(
\mathrm{CH}^{(1)}_{s}\left(
\sigma^{\mathbf{Pth}_{\boldsymbol{\mathcal{A}}}}\left(
\left(
\mathfrak{P}_{j}
\right)_{j\in\bb{\mathbf{s}}}
\right)
\right)\right)\right)
$
\allowdisplaybreaks
\begin{align*}
\qquad&=
\mathrm{pr}^{\equiv^{[1]}}_{s}\left(
\sigma^{\mathbf{Pth}_{\boldsymbol{\mathcal{A}}}}\left(
\left(
\mathrm{ip}^{(1,X)@}_{s}\left(
\mathrm{CH}^{(1)}_{s_{j}}\left(
\mathfrak{P}_{j}
\right)
\right)\right)_{j\in\bb{\mathbf{s}}}
\right)\right)
\tag{1}
\\&=
\sigma^{\mathbf{T}_{\boldsymbol{\mathcal{E}}^{\boldsymbol{\mathcal{A}}}}
(\mathbf{Pth}_{\boldsymbol{\mathcal{A}}})}\left(
\left(
\mathrm{pr}^{\equiv^{[1]}}_{s_{j}}\left(
\mathrm{ip}^{(1,X)@}_{s}\left(
\mathrm{CH}^{(1)}_{s_{j}}\left(
\mathfrak{P}_{j}
\right)
\right)
\right)
\right)_{j\in\bb{\mathbf{s}}}\right).
\tag{2}
\end{align*}
\end{flushleft}

In the just stated chain of equalities, the first equality follows from Lemma~\ref{LIpCHSigma}; finally, the last equality follows from the fact that $\mathrm{pr}^{\equiv^{[1]}}$ is a $\Sigma^{\boldsymbol{\mathcal{A}}}$-homomorphism, according to Remark~\ref{RVar}.

\textsf{Compatibility with the rewrite rules in $\mathcal{A}$.}

Let $s$ be a sort in $S$ and let $\mathfrak{p}$ be a rewrite rule in $\mathcal{A}_{s}$.

The following chain of equalities holds.
\allowdisplaybreaks
\begin{align*}
\mathrm{pr}^{\equiv^{[1]}}_{s}\left(
\mathrm{ip}^{(1,X)@}_{s}\left(
\mathrm{CH}^{(1)}_{s}\left(
\mathfrak{p}^{\mathbf{Pth}_{\boldsymbol{\mathcal{A}}}}
\right)\right)\right)&=
\mathrm{pr}^{\equiv^{[1]}}_{s}\left(
\mathfrak{p}^{\mathbf{Pth}_{\boldsymbol{\mathcal{A}}}}\right)
\tag{1}
\\&=
\mathfrak{p}^{\mathbf{T}_{\boldsymbol{\mathcal{E}}^{\boldsymbol{\mathcal{A}}}}
(\mathbf{Pth}_{\boldsymbol{\mathcal{A}}})}.
\tag{2}
\end{align*}

In the just stated chain of equalities, the first equality follows from Lemma~\ref{LIpCHEch}; finally, the last equality follows from the fact that $\mathrm{pr}^{\equiv^{[1]}}$ is a $\Sigma^{\boldsymbol{\mathcal{A}}}$-homomorphism, according to Remark~\ref{RVar}.

\textsf{Compatibility with the $0$-source.}

Let $s$ be a sort in $S$ and let $\mathfrak{P}$ be a path in $\mathrm{Pth}_{\boldsymbol{\mathcal{A}},s}$. 

The following chain of equalities holds
\begin{flushleft}
$\mathrm{pr}^{\equiv^{[1]}}_{s}\left(
\mathrm{ip}^{(1,X)@}_{s}\left(
\mathrm{CH}^{(1)}_{s}\left(
\mathrm{sc}^{0\mathbf{Pth}_{\boldsymbol{\mathcal{A}}}}_{s}\left(
\mathfrak{P}
\right)
\right)\right)\right)
$
\allowdisplaybreaks
\begin{align*}
\qquad&=
\mathrm{pr}^{\equiv^{[1]}}_{s}\left(
\mathrm{sc}^{0\mathbf{Pth}_{\boldsymbol{\mathcal{A}}}}_{s}\left(
\mathrm{ip}^{(1,X)@}_{s}\left(
\mathrm{CH}^{(1)}_{s}\left(
\mathfrak{P}
\right)
\right)\right)\right)
\tag{1}
\\&=
\mathrm{sc}^{0\mathbf{T}_{\boldsymbol{\mathcal{E}}^{\boldsymbol{\mathcal{A}}}}
(\mathbf{Pth}_{\boldsymbol{\mathcal{A}}})}_{s}\left(
\mathrm{pr}^{\equiv^{[1]}}_{s}\left(
\mathrm{ip}^{(1,X)@}_{s}\left(
\mathrm{CH}^{(1)}_{s}\left(
\mathfrak{P}
\right)
\right)\right)\right).
\tag{2}
\end{align*}
\end{flushleft}

In the just stated chain of equalities, the first equality follows from Lemma~\ref{LIpCHScTg}; finally, the last equality follows from the fact that $\mathrm{pr}^{\equiv^{[1]}}$ is a $\Sigma^{\boldsymbol{\mathcal{A}}}$-homomorphism, according to Remark~\ref{RVar}.

\textsf{Compatibility with the $0$-target.}

Let $s$ be a sort in $S$ and let $\mathfrak{P}$ be a path in $\mathrm{Pth}_{\boldsymbol{\mathcal{A}},s}$.  

The following equality holds
\allowdisplaybreaks
\begin{multline*}
\mathrm{pr}^{\equiv^{[1]}}_{s}\left(
\mathrm{ip}^{(1,X)@}_{s}\left(
\mathrm{CH}^{(1)}_{s}\left(
\mathrm{tg}^{0\mathbf{Pth}_{\boldsymbol{\mathcal{A}}}}_{s}\left(
\mathfrak{P}
\right)
\right)\right)\right)
\\=
\mathrm{tg}^{0\mathbf{T}_{\boldsymbol{\mathcal{E}}^{\boldsymbol{\mathcal{A}}}}
(\mathbf{Pth}_{\boldsymbol{\mathcal{A}}})}_{s}\left(
\mathrm{pr}^{\equiv^{[1]}}_{s}\left(
\mathrm{ip}^{(1,X)@}_{s}\left(
\mathrm{CH}^{(1)}_{s}\left(
\mathfrak{P}
\right)
\right)\right)\right).
\end{multline*}

The proof of this case is similar to that presented for the $0$-source.

\textsf{Compatibility with the $0$-composition.}

Let $s$ be a sort in $S$ and let $\mathfrak{P},\mathfrak{Q}$ be  paths in $\mathrm{Pth}_{\boldsymbol{\mathcal{A}},s}$ satisfying  
\[
\mathrm{sc}^{(0,1)}_{s}\left(
\mathfrak{Q}
\right)
=
\mathrm{tg}^{(0,1)}_{s}\left(
\mathfrak{P}
\right).
\]

We want to check that 
\allowdisplaybreaks
\begin{multline*}
\mathrm{pr}^{\equiv^{[1]}}_{s}\left(
\mathrm{ip}^{(1,X)@}_{s}\left(
\mathrm{CH}^{(1)}_{s}\left(
\mathfrak{Q}
\circ^{0\mathbf{Pth}_{\boldsymbol{\mathcal{A}}}}_{s}
\mathfrak{P}
\right)\right)\right)
\\=
\mathrm{pr}^{\equiv^{[1]}}_{s}\left(
\mathrm{ip}^{(1,X)@}_{s}\left(
\mathrm{CH}^{(1)}_{s}\left(
\mathfrak{Q}
\right)
\right)\right)
\circ^{0\mathbf{T}_{\boldsymbol{\mathcal{E}}^{\boldsymbol{\mathcal{A}}}}
(\mathbf{Pth}_{\boldsymbol{\mathcal{A}}})}_{s}
\mathrm{pr}^{\equiv^{[1]}}_{s}\left(
\mathrm{ip}^{(1,X)@}_{s}\left(
\mathrm{CH}^{(1)}_{s}\left(
\mathfrak{P}
\right)
\right)\right).
\end{multline*}

Developing both sides of the above equation, and taking into account Remark~\ref{RVar}, we have that to check that the following equality between classes holds
\allowdisplaybreaks
\begin{multline*}
\left[
\mathrm{ip}^{(1,X)@}_{s}\left(
\mathrm{CH}^{(1)}_{s}\left(
\mathfrak{Q}
\circ^{0\mathbf{Pth}_{\boldsymbol{\mathcal{A}}}}_{s}
\mathfrak{P}
\right)\right)
\right]_{\equiv^{[1]}}
\\=
\left[
\mathrm{ip}^{(1,X)@}_{s}\left(
\mathrm{CH}^{(1)}_{s}\left(
\mathfrak{Q}
\right)\right)
\circ^{0\mathbf{Pth}_{\boldsymbol{\mathcal{A}}}}_{s}
\mathrm{ip}^{(1,X)@}_{s}\left(
\mathrm{CH}^{(1)}_{s}\left(
\mathfrak{P}
\right)\right)
\right]_{\equiv^{[1]}}.
\end{multline*}

But this is equivalent to prove that, for every many-sorted partial $\Sigma^{\boldsymbol{\mathcal{A}}}$-algebra $\mathbf{B}$ in $\mathbf{PAlg}(\boldsymbol{\mathcal{E}}^{\boldsymbol{\mathcal{A}}})$ and every $\Sigma^{\boldsymbol{\mathcal{A}}}$-homomorphism $f\colon\mathbf{Pth}_{\boldsymbol{\mathcal{A}}}\mor \mathbf{B}$ we have that 
\allowdisplaybreaks
\begin{multline*}
f^{\mathrm{Sch}}_{s}\left(
\mathrm{ip}^{(1,X)@}_{s}\left(
\mathrm{CH}^{(1)}_{s}\left(
\mathfrak{Q}
\circ^{0\mathbf{Pth}_{\boldsymbol{\mathcal{A}}}}_{s}
\mathfrak{P}
\right)\right)
\right)
\\=
f^{\mathrm{Sch}}_{s}\left(
\mathrm{ip}^{(1,X)@}_{s}\left(
\mathrm{CH}^{(1)}_{s}\left(
\mathfrak{Q}
\right)\right)
\circ^{0\mathbf{Pth}_{\boldsymbol{\mathcal{A}}}}_{s}
\mathrm{ip}^{(1,X)@}_{s}\left(
\mathrm{CH}^{(1)}_{s}\left(
\mathfrak{P}
\right)\right)
\right).
\end{multline*}

Let us recall that, by Proposition~\ref{PIpCH}, we are dealing with  proper paths. Hence, the last mentioned equality reduces to
\allowdisplaybreaks
\begin{multline*}
f_{s}\left(
\mathrm{ip}^{(1,X)@}_{s}\left(
\mathrm{CH}^{(1)}_{s}\left(
\mathfrak{Q}
\circ^{0\mathbf{Pth}_{\boldsymbol{\mathcal{A}}}}_{s}
\mathfrak{P}
\right)\right)
\right)
\\=
f_{s}\left(
\mathrm{ip}^{(1,X)@}_{s}\left(
\mathrm{CH}^{(1)}_{s}\left(
\mathfrak{Q}
\right)\right)
\circ^{0\mathbf{Pth}_{\boldsymbol{\mathcal{A}}}}_{s}
\mathrm{ip}^{(1,X)@}_{s}\left(
\mathrm{CH}^{(1)}_{s}\left(
\mathfrak{P}
\right)\right)
\right).
\tag{E2}\label{PVarKerE2}
\end{multline*}

We provide the diagram in Figure~\ref{FTheBigOneP} to better understand the elements under consideration.

\begin{figure}
\begin{center}
\begin{tikzpicture}
[ACliment/.style={-{To [angle'=45, length=5.75pt, width=4pt, round]}
}, scale=0.8]
\node[] (X) 		at 	(-3,6) 	[] 	{$\mathbf{D}_{\Sigma^{\boldsymbol{\mathcal{A}}}}(X)$};
\node[] (Sdg) 	at 	(0,5) 	[] 	{$\mathbf{Sch}(\mathrm{ip}^{(1,X)})$};
\node[]	(T2)		at 	(-3,3)	[]	{$\mathbf{Pth}_{\boldsymbol{\mathcal{A}}}$};
\node[]	(Sf)		at 	(0,0)	[]	{$\mathbf{Sch}(f)$};
\node[]	(T)		at 	(8,2)	[]	{$\mathbf{T}_{\Sigma^{\boldsymbol{\mathcal{A}}}}(X)$};
\node[]	(F2)		at 	(5,0)	[]	{$\mathbf{F}_{\Sigma^{\boldsymbol{\mathcal{A}}}}(\mathbf{Pth}_{\boldsymbol{\mathcal{A}}})$};
\node[]	(PK)		at 	(2.5,2.5)	[]	{$[\mathbf{Pth}_{\boldsymbol{\mathcal{A}}}]$};
\node[]	(FB)		at 	(5,-2)	[]	{$\mathbf{F}_{\Sigma^{\boldsymbol{\mathcal{A}}}}(\mathbf{B})$};
\node[]	(B)		at 	(0,-2)	[]	{$\mathbf{B}$};
\draw[ACliment, bend right]  (X) 	to node [below left]	
{$\mathrm{ip}^{(1,X)}$} (T2);
\draw[ACliment, bend left]  (X) 	to node [above right]	{$\eta^{(1,X)}$} (T);
\draw[ACliment]  (X) 	to node [above right]	{$\mathrm{in}^{X}$} (Sdg);
\draw[ACliment, bend right]  (T2) 	to node [below left]	{$f$} (B);
\draw[ACliment]  (T2) 	to node [midway, fill=white]	{$\mathrm{in}^{\mathbf{Pth}_{\boldsymbol{\mathcal{A}}}}$} (Sf);
\draw[ACliment]  (Sf) 	to node [above]	{$\mathrm{in}^{\mathbf{Sch}(f)}$} (F2);
\draw[ACliment]  (Sdg) 	to node [midway, fill=white]	{$\mathrm{in}^{\mathbf{Sch}(\mathrm{ip}^{(1,X)})}$} (T);
\draw[ACliment]  (T) 	to node [below right]	 {$\mathrm{ip}^{(1,X)@}$} (F2);
\draw[ACliment]  (Sdg) 	to node [midway, fill=white]	{$\mathrm{ip}^{(1,X)\mathrm{Sch}}$} (T2);
\draw[ACliment]  (B) 	to node [below]	{$\eta^{(1,\mathbf{B})}$} (FB);
\draw[ACliment] (Sf) to node [right] {$f^{\mathrm{Sch}}$} (B);
\draw[ACliment] (F2) to node [right] {$f^{@}$} (FB);
\draw[ACliment] (PK) to node [below] {$\mathrm{CH}^{(1)\mathrm{m}}$} (T);
\draw[ACliment] (T2) to node [above] {$\mathrm{pr}^{\mathrm{Ker}(\mathrm{CH}^{(1)})}$} (PK);
\draw[ACliment]  (T2) 	to node [midway, fill=white]	{$\eta^{(1,\mathbf{Pth}_{\boldsymbol{\mathcal{A}}})}$} (F2);
\end{tikzpicture}
\end{center}
\caption{The first big diagram for paths.}
\label{FTheBigOneP}
\end{figure}

Before presenting a proof for the general case, we will see a proof of Equation~\ref{PVarKerE2} when one of the  paths involved in the $0$-composition is a $(1,0)$-identity path.

\begin{claim}\label{CVarKer} Let $\mathfrak{P}$ and $\mathfrak{Q}$ be two paths in $\mathrm{Pth}_{\boldsymbol{\mathcal{A}},s}$ such that $\mathrm{sc}^{(0,1)}_{s}(\mathfrak{Q})=\mathrm{tg}^{(0,1)}_{s}(\mathfrak{P})$. If either $\mathfrak{P}$ or $\mathfrak{Q}$ is a $(1,0)$-identity  path, then for every many-sorted partial $\Sigma^{\boldsymbol{\mathcal{A}}}$-algebra $\mathbf{B}$ in $\mathbf{PAlg}(\boldsymbol{\mathcal{E}}^{\boldsymbol{\mathcal{A}}})$ and every $\Sigma^{\boldsymbol{\mathcal{A}}}$-homomorphism $f\colon \mathbf{Pth}_{\boldsymbol{\mathcal{A}}}\mor \mathbf{B}$, Equation~\ref{PVarKerE2} holds.
\end{claim}

Assume without loss of generality that $\mathfrak{P}$ is a $(1,0)$-identity  path, then the following chain of equalities holds.
\begin{flushleft}
$f_{s}\left(
\mathrm{ip}^{(1,X)@}_{s}\left(
\mathrm{CH}^{(1)}_{s}\left(
\mathfrak{Q}
\circ^{0\mathbf{Pth}_{\boldsymbol{\mathcal{A}}}}_{s}
\mathfrak{P}
\right)\right)
\right)$
\allowdisplaybreaks
\begin{align*}
&=f_{s}\left(
\mathrm{ip}^{(1,X)@}_{s}\left(
\mathrm{CH}^{(1)}_{s}\left(
\mathfrak{Q}
\right)\right)
\right)
\tag{1}
\\&=f_{s}\left(
\mathrm{ip}^{(1,X)@}_{s}\left(
\mathrm{CH}^{(1)}_{s}\left(
\mathfrak{Q}
\right)\right)
\circ^{0\mathbf{Pth}_{\boldsymbol{\mathcal{A}}}}_{s}
\mathfrak{P}
\right)
\tag{2}
\\&=f_{s}\left(
\mathrm{ip}^{(1,X)@}_{s}\left(
\mathrm{CH}^{(1)}_{s}\left(
\mathfrak{Q}
\right)\right)
\right)
\circ^{0\mathbf{B}}_{s}
f_{s}\left(
\mathfrak{P}
\right)
\tag{3}
\\&=
f_{s}\left(
\mathrm{ip}^{(1,X)@}_{s}\left(
\mathrm{CH}^{(1)}_{s}\left(
\mathfrak{Q}
\right)\right)
\right)
\circ^{0\mathbf{B}}_{s}
f_{s}\left(
\mathrm{ip}^{(1,X)@}_{s}\left(
\mathrm{CH}^{(1)}_{s}\left(
\mathfrak{P}
\right)\right)
\right).
\tag{4}
\end{align*}
\end{flushleft}

Note that since $\mathfrak{P}$ is a $(1,0)$-identity path, this path must be equal to the $(1,0)$-identity  path on the $(0,1)$-source of $\mathfrak{Q}$, so the $0$-composition of the paths is well-defined. This justifies the first equality. The second equality follows from the fact that, according to Proposition~\ref{PIpCH}, $\mathrm{ip}^{(1,X)@}_{s}(\mathrm{CH}^{(1)}_{s}(\mathfrak{Q}))$ is a path in $[\mathfrak{Q}]_{s}$. Thus, by Lemma~\ref{LCH}, the  paths $\mathrm{ip}^{(1,X)@}_{s}(\mathrm{CH}^{(1)}_{s}(\mathfrak{Q}))$ and $\mathfrak{Q}$ have the same $(0,1)$-source. Moreover, since $\mathfrak{P}$ is the  $(1,0)$-identity  path on the $(0,1)$-source of $\mathfrak{Q}$, we have that
\[
\mathrm{ip}^{(1,X)@}_{s}\left(
\mathrm{CH}^{(1)}_{s}\left(
\mathfrak{Q}
\right)\right)
=
\mathrm{ip}^{(1,X)@}_{s}\left(
\mathrm{CH}^{(1)}_{s}\left(
\mathfrak{Q}
\right)\right)
\circ^{0\mathbf{Pth}_{\boldsymbol{\mathcal{A}}}}_{s}
\mathfrak{P};
\]
The third equality follows from the fact that, by assumption, $f$ is a $\Sigma^{\boldsymbol{\mathcal{A}}}$-homomorphism; finally, the last equality follows from Corollary~\ref{CIpId}. In this regard, since $\mathfrak{P}$ is a $(1,0)$-identity  path, we have that $\mathrm{ip}^{(1,X)@}_{s}(
\mathrm{CH}^{(1)}_{s}(
\mathfrak{P}
))=\mathfrak{P}$.

The same argument applies in case $\mathfrak{Q}$ is a $(1,0)$-identity  path. In this case, we will argue taking into account that $\mathfrak{Q}$ must be the $(1,0)$-identity path on the $(0,1)$-target of $\mathfrak{P}$.

This proves Claim~\ref{CVarKer}.

We are now in position to prove the general case on $\mathfrak{Q}\circ^{0\mathbf{Pth}_{\boldsymbol{\mathcal{A}}}}_{s}\mathfrak{P}$ by Artinian induction on $(\coprod\mathrm{Pth}_{\boldsymbol{\mathcal{A}}}, \leq_{\mathbf{Pth}_{\boldsymbol{\mathcal{A}}}})$.

\textsf{Base step of the Artinian induction.}

Let $(\mathfrak{Q}\circ^{0\mathbf{Pth}_{\boldsymbol{\mathcal{A}}}}_{s}\mathfrak{P}, s)$ be a minimal element in $(\coprod\mathrm{Pth}_{\boldsymbol{\mathcal{A}}}, \leq_{\mathbf{Pth}_{\boldsymbol{\mathcal{A}}}})$. Then by Proposition~\ref{PMinimal}, the path $\mathfrak{Q}\circ^{0\mathbf{Pth}_{\boldsymbol{\mathcal{A}}}}_{s}\mathfrak{P}$ is either a $(1,0)$-identity path or an echelon. In any case, either $\mathfrak{P}$ or $\mathfrak{Q}$ must be a $(1,0)$-identity  path. The statement follows by Claim~\ref{CVarKer}.

\textsf{Inductive step of the Artinian induction.}

Let $(\mathfrak{Q}\circ^{0\mathbf{Pth}_{\boldsymbol{\mathcal{A}}}}_{s}\mathfrak{P},s)$ be a non-minimal element in $(\coprod\mathrm{Pth}_{\boldsymbol{\mathcal{A}}}, \leq_{\mathbf{Pth}_{\boldsymbol{\mathcal{A}}}})$. Let us suppose that, for every sort  $t\in S$ and every path $\mathfrak{Q}'\circ^{0\mathbf{Pth}_{\boldsymbol{\mathcal{A}}}}_{t}\mathfrak{P}'$ in $\mathrm{Pth}_{\boldsymbol{\mathcal{A}},t}$, if $(\mathfrak{Q}'\circ^{0\mathbf{Pth}_{\boldsymbol{\mathcal{A}}}}_{t}\mathfrak{P}',t)$ $<_{\mathbf{Pth}_{\boldsymbol{\mathcal{A}}}}$-precedes $(\mathfrak{Q}\circ^{0\mathbf{Pth}_{\boldsymbol{\mathcal{A}}}}_{s}\mathfrak{P},s)$, then the statement holds for $\mathfrak{Q}'\circ^{0\mathbf{Pth}_{\boldsymbol{\mathcal{A}}}}_{t}\mathfrak{P}'$, i.e., the following equality holds
\allowdisplaybreaks
\begin{multline*}
f_{t}\left(
\mathrm{ip}^{(1,X)@}_{t}\left(
\mathrm{CH}^{(1)}_{t}\left(
\mathfrak{Q}'
\circ^{0\mathbf{Pth}_{\boldsymbol{\mathcal{A}}}}_{t}
\mathfrak{P}'
\right)\right)\right)
\\=
f_{t}\left(
\mathrm{ip}^{(1,X)@}_{t}\left(
\mathrm{CH}^{(1)}_{t}\left(
\mathfrak{Q}'
\right)\right)\right)
\circ^{0\mathbf{B}}_{t}
f_{t}\left(
\mathrm{ip}^{(1,X)@}_{t}\left(
\mathrm{CH}^{(1)}_{t}\left(
\mathfrak{P}'
\right)\right)\right).
\end{multline*}

Since $(\mathfrak{Q}\circ^{0\mathbf{Pth}_{\boldsymbol{\mathcal{A}}}}_{s}\mathfrak{P},s)$ is a non-minimal element in $(\coprod\mathrm{Pth}_{\boldsymbol{\mathcal{A}}}, \leq_{\mathbf{Pth}_{\boldsymbol{\mathcal{A}}}})$ and, by Claim~\ref{CVarKer}, we can assume that neither $\mathfrak{Q}$ nor $\mathfrak{P}$ are $(1,0)$-identity paths, we have, by Lemma~\ref{LOrdI}, that $\mathfrak{Q}\circ^{0\mathbf{Pth}_{\boldsymbol{\mathcal{A}}}}_{s}\mathfrak{P}$  it is either~(1) a path of length strictly greater than one containing at least one echelon or~(2) an echelonless path.

If~(1), then let $i\in \bb{\mathfrak{Q}\circ^{0\mathbf{Pth}_{\boldsymbol{\mathcal{A}}}}_{s}\mathfrak{P}}$ be the first index for which the one-step subpath $(\mathfrak{Q}\circ^{0\mathbf{Pth}_{\boldsymbol{\mathcal{A}}}}_{s}\mathfrak{P})^{i,i}$ of $\mathfrak{Q}\circ^{0\mathbf{Pth}_{\boldsymbol{\mathcal{A}}}}_{s}\mathfrak{P}$ is an echelon.  We distinguish the cases~(1.1) $i=0$ and~(1.2) $i>0$.

If~(1.1), i.e., if $i=0$, since we are assuming that $\mathfrak{P}$ is not a $(1,0)$-identity path, we have that $\mathfrak{P}$ has an echelon in its first step. Then it could be the case that either (1.1.1) $\mathfrak{P}$ is an echelon or (1.1.2) $\mathfrak{P}$ is a path of length strictly greater than one containing an echelon on its first step.

If (1.1.1), i.e., if we consider the case in which $\mathfrak{P}$ is an echelon, then the value of the Curry-Howard mapping at $\mathfrak{Q}\circ^{0\mathbf{Pth}_{\boldsymbol{\mathcal{A}}}}_{s}\mathfrak{P}$ is given by
$$
\mathrm{CH}^{(1)}_{s}\left(
\mathfrak{Q}\circ^{0\mathbf{Pth}_{\boldsymbol{\mathcal{A}}}}_{s}\mathfrak{P}
\right)=
\mathrm{CH}^{(1)}_{s}\left(
\mathfrak{Q}
\right)
\circ_{s}^{0\mathbf{T}_{\Sigma^{\boldsymbol{\mathcal{A}}}}(X)}
\mathrm{CH}^{(1)}_{s}\left(
\mathfrak{P}
\right).
$$

Therefore, the following chain of equalities holds
\begin{flushleft}
$
f_{s}\left(
\mathrm{ip}^{(1,X)@}_{s}\left(
\mathrm{CH}^{(1)}_{s}\left(
\mathfrak{Q}
\circ^{0\mathbf{Pth}_{\boldsymbol{\mathcal{A}}}}_{s}
\mathfrak{P}
\right)\right)\right)=
$
\allowdisplaybreaks
\begin{align*}
\qquad&=
f_{s}\left(
\mathrm{ip}^{(1,X)@}_{s}\left(
\mathrm{CH}^{(1)}_{s}\left(
\mathfrak{Q}
\right)
\circ^{0\mathbf{T}_{\Sigma^{\boldsymbol{\mathcal{A}}}}(X)}_{s}
\mathrm{CH}^{(1)}_{s}\left(
\mathfrak{P}
\right)\right)\right)
\tag{1}
\\&=
f_{s}\left(
\mathrm{ip}^{(1,X)@}_{s}\left(
\mathrm{CH}^{(1)}_{s}\left(
\mathfrak{Q}
\right)\right)
\circ_{s}^{0\mathbf{F}_{\Sigma^{\boldsymbol{\mathcal{A}}}}
(\mathbf{Pth}_{\boldsymbol{\mathcal{A}}})}
\mathrm{ip}^{(1,X)@}_{s}\left(
\mathrm{CH}^{(1)}_{s}\left(
\mathfrak{P}
\right)\right)\right)
\tag{2}
\\&=
f_{s}\left(
\mathrm{ip}^{(1,X)@}_{s}\left(
\mathrm{CH}^{(1)}_{s}\left(
\mathfrak{Q}
\right)\right)
\circ_{s}^{0\mathbf{Pth}_{\boldsymbol{\mathcal{A}}}}
\mathrm{ip}^{(1,X)@}_{s}\left(
\mathrm{CH}^{(1)}_{s}\left(
\mathfrak{P}
\right)\right)\right)
\tag{3}
\\&=
f_{s}\left(
\mathrm{ip}^{(1,X)@}_{s}\left(
\mathrm{CH}^{(1)}_{s}\left(
\mathfrak{Q}
\right)\right)\right)
\circ_{s}^{0\mathbf{B}}
f_{s}\left(
\mathrm{ip}^{(1,X)@}_{s}\left(
\mathrm{CH}^{(1)}_{s}\left(
\mathfrak{P}
\right)\right)\right).
\tag{4}
\end{align*}
\end{flushleft}

The first equality follows from the previous discussion on the value of the Curry-Howard mapping at $\mathfrak{Q}\circ^{0\mathbf{Pth}_{\boldsymbol{\mathcal{A}}}}_{s}\mathfrak{P}$; the second equality follows from the fact that, by Definition~\ref{DIp}, $\mathrm{ip}^{(1,X)@}$ is a $\Sigma^{\boldsymbol{\mathcal{A}}}$-homomorphism; the third equality follows from the fact that, in virtue of Proposition~\ref{PIpCH}, both $\mathrm{ip}^{(1,X)@}_{s}(
\mathrm{CH}^{(1)}_{s}(\mathfrak{Q}))$ and $\mathrm{ip}^{(1,X)@}_{s}(
\mathrm{CH}^{(1)}_{s}(\mathfrak{P}))$ are paths. Thus, the interpretation of the operation symbol for the $0$-composition, $\circ^{0}_{s}$, in the $\Sigma^{\boldsymbol{\mathcal{A}}}$-algebra $\mathbf{F}_{\Sigma^{\boldsymbol{\mathcal{A}}}}(\mathbf{Pth}_{\boldsymbol{\mathcal{A}}})$ becomes the corresponding interpretation of the $0$-composition $\circ^{0}_{s}$ in the $\Sigma^{\boldsymbol{\mathcal{A}}}$-algebra $\mathbf{Pth}_{\boldsymbol{\mathcal{A}}}$;  finally, the last equality follows from the fact that, by assumption, $f$ is a $\Sigma^{\boldsymbol{\mathcal{A}}}$-homomorphism.

The case (1.1.1) of $\mathfrak{P}$ being an echelon follows.

If~(1.1.2), i.e., if we are in the case in which $\mathfrak{P}$ is a path of length strictly greater than one containing an echelon on its first step, then the value of the Curry-Howard mapping at $\mathfrak{P}$ is given by
$$
\mathrm{CH}^{(1)}_{s}\left(
\mathfrak{P}
\right)=
\mathrm{CH}^{(1)}_{s}
\left(
\mathfrak{P}^{1,\bb{\mathfrak{P}}-1}
\right)
\circ_{s}^{0\mathbf{T}_{\Sigma^{\boldsymbol{\mathcal{A}}}}(X)}
\mathrm{CH}^{(1)}_{s}\left(
\mathfrak{P}^{0,0}
\right).
$$

Moreover, $(\mathfrak{Q}\circ^{0\mathbf{Pth}_{\boldsymbol{\mathcal{A}}}}_{s}\mathfrak{P})^{1,\bb{\mathfrak{Q}\circ^{0\mathbf{Pth}_{\boldsymbol{\mathcal{A}}}}_{s}\mathfrak{P}}-1}=\mathfrak{Q}\circ_{s}^{0\mathbf{Pth}_{\boldsymbol{\mathcal{A}}}}\mathfrak{P}^{1,\bb{\mathfrak{P}}-1}$ Hence, the value of the Curry-Howard mapping at $\mathfrak{Q}\circ^{0\mathbf{Pth}_{\boldsymbol{\mathcal{A}}}}_{s}\mathfrak{P}$ is given by
$$
\mathrm{CH}^{(1)}_{s}
\left(
\mathfrak{Q}\circ^{0\mathbf{Pth}_{\boldsymbol{\mathcal{A}}}}_{s}\mathfrak{P}
\right)=
\mathrm{CH}^{(1)}_{s}\left(
\mathfrak{Q}\circ^{0\mathbf{Pth}_{\boldsymbol{\mathcal{A}}}}_{s}\mathfrak{P}
^{1,\bb{\mathfrak{P}}-1}\right)
\circ_{s}^{0\mathbf{T}_{\Sigma^{\boldsymbol{\mathcal{A}}}}(X)}
\mathrm{CH}^{(1)}_{s}\left(\mathfrak{P}
^{0,0}\right).
$$

Therefore, the following chain of equalities holds
\begin{flushleft}
$
f_{s}\left(
\mathrm{ip}^{(1,X)@}_{s}\left(
\mathrm{CH}^{(1)}_{s}\left(
\mathfrak{Q}
\circ^{0\mathbf{Pth}_{\boldsymbol{\mathcal{A}}}}_{s}
\mathfrak{P}
\right)\right)\right)
$
\allowdisplaybreaks
\begin{align*}
\qquad&=
f_{s}\left(
\mathrm{ip}^{(1,X)@}_{s}\left(
\mathrm{CH}^{(1)}_{s}\left(
\mathfrak{Q}\circ^{0\mathbf{Pth}_{\boldsymbol{\mathcal{A}}}}_{s}\mathfrak{P}
^{1,\bb{\mathfrak{P}}-1}
\right)
\circ_{s}^{0\mathbf{T}_{\Sigma^{\boldsymbol{\mathcal{A}}}}(X)}
\mathrm{CH}^{(1)}_{s}\left(
\mathfrak{P}
^{0,0}
\right)\right)\right)
\tag{1}
\\&=
f_{s}\Big(
\mathrm{ip}^{(1,X)@}_{s}\left(
\mathrm{CH}^{(1)}_{s}\left(
\mathfrak{Q}\circ^{0\mathbf{Pth}_{\boldsymbol{\mathcal{A}}}}_{s}\mathfrak{P}
^{1,\bb{\mathfrak{P}}-1}\right)\right)
\circ_{s}^{0\mathbf{F}_{\Sigma^{\boldsymbol{\mathcal{A}}}}
(\mathbf{Pth}_{\boldsymbol{\mathcal{A}}})}
\\&\qquad\qquad\qquad\qquad\qquad\qquad\qquad\qquad\qquad\qquad
\mathrm{ip}^{(1,X)@}_{s}\left(
\mathrm{CH}^{(1)}_{s}\left(
\mathfrak{P}^{0,0}
\right)\right)\Big)
\tag{2}
\\&=
f_{s}\left(
\mathrm{ip}^{(1,X)@}_{s}\left(
\mathrm{CH}^{(1)}_{s}\left(
\mathfrak{Q}\circ^{0\mathbf{Pth}_{\boldsymbol{\mathcal{A}}}}_{s}\mathfrak{P}
^{1,\bb{\mathfrak{P}}-1}
\right)\right)
\circ_{s}^{0\mathbf{Pth}_{\boldsymbol{\mathcal{A}}}}
\right.
\\&\qquad\qquad\qquad\qquad\qquad\qquad\qquad\qquad\qquad\qquad
\left.
\mathrm{ip}^{(1,X)@}_{s}\left(
\mathrm{CH}^{(1)}_{s}\left(
\mathfrak{P}^{0,0}\right)\right)\right)
\tag{3}
\\&=
f_{s}\left(
\mathrm{ip}^{(1,X)@}_{s}\left(
\mathrm{CH}^{(1)}_{s}\left(
\mathfrak{Q}\circ^{0\mathbf{Pth}_{\boldsymbol{\mathcal{A}}}}_{s}\mathfrak{P}
^{1,\bb{\mathfrak{P}}-1}
\right)\right)\right)
\circ_{s}^{0\mathbf{B}}
\\&\qquad\qquad\qquad\qquad\qquad\qquad\qquad\qquad\qquad
f_{s}\left(
\mathrm{ip}^{(1,X)@}_{s}\left(
\mathrm{CH}^{(1)}_{s}\left(
\mathfrak{P}^{0,0}
\right)\right)\right)
\tag{4}
\\&=
\left(
f_{s}\left(
\mathrm{ip}^{(1,X)@}_{s}\left(
\mathrm{CH}^{(1)}_{s}\left(
\mathfrak{Q}
\right)\right)\right)
\circ_{s}^{0\mathbf{B}}
f_{s}\left(
\mathrm{ip}^{(1,X)@}_{s}\left(
\mathrm{CH}^{(1)}_{s}\left(
\mathfrak{P}^{1,\bb{\mathfrak{P}}-1}
\right)\right)\right)
\right)
\circ_{s}^{0\mathbf{B}}
\\&\qquad\qquad\qquad\qquad\qquad\qquad\qquad\qquad\qquad
f_{s}\left(
\mathrm{ip}^{(1,X)@}_{s}\left(
\mathrm{CH}^{(1)}_{s}\left(
\mathfrak{P}^{0,0}
\right)\right)\right)
\tag{5}
\\&=
f_{s}\left(
\mathrm{ip}^{(1,X)@}_{s}\left(
\mathrm{CH}^{(1)}_{s}\left(
\mathfrak{Q}
\right)\right)\right)
\circ_{s}^{0\mathbf{B}}
\left(
f_{s}\left(
\mathrm{ip}^{(1,X)@}_{s}\left(
\mathrm{CH}^{(1)}_{s}\left(
\mathfrak{P}^{1,\bb{\mathfrak{P}}-1}
\right)\right)\right)
\circ_{s}^{0\mathbf{B}}
\right.
\\&\qquad\qquad\qquad\qquad\qquad\qquad\qquad\qquad
\left.
f_{s}\left(
\mathrm{ip}^{(1,X)@}_{s}\left(
\mathrm{CH}^{(1)}_{s}\left(
\mathfrak{P}^{0,0}
\right)\right)\right)
\right)
\tag{6}
\\&=
f_{s}\left(
\mathrm{ip}^{(1,X)@}_{s}\left(
\mathrm{CH}^{(1)}_{s}\left(
\mathfrak{Q}
\right)\right)\right)
\circ_{s}^{0\mathbf{B}}
f_{s}
\left(
\mathrm{ip}^{(1,X)@}_{s}\left(
\mathrm{CH}^{(1)}_{s}\left(
\mathfrak{P}^{1,\bb{\mathfrak{P}}-1}
\right)\right)
\circ_{s}^{0\mathbf{Pth}_{\boldsymbol{\mathcal{A}}}}
\right.
\\&\qquad\qquad\qquad\qquad\qquad\qquad\qquad\qquad\qquad\qquad
\left.
\mathrm{ip}^{(1,X)@}_{s}\left(
\mathrm{CH}^{(1)}_{s}\left(\mathfrak{P}^{0,0}
\right)\right)
\right)
\tag{7}
\\&=
f_{s}\left(
\mathrm{ip}^{(1,X)@}_{s}\left(
\mathrm{CH}^{(1)}_{s}\left(
\mathfrak{Q}
\right)\right)\right)
\circ_{s}^{0\mathbf{B}}
f_{s}
\left(
\mathrm{ip}^{(1,X)@}_{s}\left(
\mathrm{CH}^{(1)}_{s}\left(
\mathfrak{P}^{1,\bb{\mathfrak{P}}-1}
\right)\right)
\right.
\\&\qquad\qquad\qquad\qquad\qquad\qquad\qquad
\left.
\circ_{s}^{0\mathbf{F}_{\Sigma^{\boldsymbol{\mathcal{A}}}}(
\mathbf{Pth}_{\boldsymbol{\mathcal{A}}})}
\mathrm{ip}^{(1,X)@}_{s}\left(
\mathrm{CH}^{(1)}_{s}\left(
\mathfrak{P}^{0,0}
\right)\right)
\right)
\tag{8}
\\&=
f_{s}\left(
\mathrm{ip}^{(1,X)@}_{s}\left(
\mathrm{CH}^{(1)}_{s}\left(
\mathfrak{Q}
\right)\right)\right)
\circ_{s}^{0\mathbf{B}}
f_{s}\left(
\mathrm{ip}^{(1,X)@}_{s}\left(
\mathrm{CH}^{(1)}_{s}\left(
\mathfrak{P}^{1,\bb{\mathfrak{P}}-1}
\right)
\right.
\right.
\\&\qquad\qquad\qquad\qquad\qquad\qquad\qquad\qquad\qquad\qquad
\left.
\left.
\circ_{s}^{0\mathbf{T}_{\Sigma^{\boldsymbol{\mathcal{A}}}}(X)}
\mathrm{CH}^{(1)}_{s}\left(
\mathfrak{P}^{0,0}
\right)\right)\right)
\tag{9}
\\&=
f_{s}\left(
\mathrm{ip}^{(1,X)@}_{s}\left(
\mathrm{CH}^{(1)}_{s}\left(
\mathfrak{Q}
\right)\right)\right)
\circ_{s}^{0\mathbf{B}}
f_{s}\left(
\mathrm{ip}^{(1,X)@}_{s}\left(
\mathrm{CH}^{(1)}_{s}\left(
\mathfrak{P}
\right)\right)\right).
\tag{10}
\end{align*}
\end{flushleft}

The first equality follows from the previous discussion on the value of the Curry-Howard mapping at $\mathfrak{Q}\circ^{0\mathbf{Pth}_{\boldsymbol{\mathcal{A}}}}_{s}\mathfrak{P}$; the second equality follows from the fact that, by Definition~\ref{DIp}, $\mathrm{ip}^{(1,X)@}$ is a $\Sigma^{\boldsymbol{\mathcal{A}}}$-homomorphism; the third equality follows from the fact that, by Proposition~\ref{PIpCH}, both $\mathrm{ip}^{(1,X)@}_{s}(
\mathrm{CH}^{(1)}_{s}(\mathfrak{Q}
\circ_{s}^{0\mathbf{Pth}_{\boldsymbol{\mathcal{A}}}}\mathfrak{P}^{1,\bb{\mathfrak{P}}-1}
))$ and $\mathrm{ip}^{(1,X)@}_{s}(
\mathrm{CH}^{(1)}_{s}(\mathfrak{P}^{0,0}))$ are paths. Thus, the interpretation of the operation symbol for the $0$-composition, $\circ^{0}_{s}$,  in the $\Sigma^{\boldsymbol{\mathcal{A}}}$-algebra $\mathbf{F}_{\Sigma^{\boldsymbol{\mathcal{A}}}}(\mathbf{Pth}_{\boldsymbol{\mathcal{A}}})$ becomes the corresponding interpretation of the $0$-composition $\circ^{0}_{s}$ in the $\Sigma^{\boldsymbol{\mathcal{A}}}$-algebra $\mathbf{Pth}_{\boldsymbol{\mathcal{A}}}$; the fourth equality follows from the fact that, by assumption, $f$ is a $\Sigma^{\boldsymbol{\mathcal{A}}}$-homomorphism; the fifth equality follows by induction. Let us note that the pair $(\mathfrak{Q}\circ^{0\mathbf{Pth}_{\boldsymbol{\mathcal{A}}}}_{s}\mathfrak{P}^{1,\bb{\mathfrak{P}}-1},s)$ $\prec_{\mathbf{Pth}_{\boldsymbol{\mathcal{A}}}}$-precedes $(\mathfrak{Q}\circ^{0\mathbf{Pth}_{\boldsymbol{\mathcal{A}}}}_{s}\mathfrak{P},s)$. Hence we have that 
\begin{multline*}
f_{s}\left(
\mathrm{ip}^{(1,X)@}_{s}\left(
\mathrm{CH}^{(1)}_{s}\left(
\mathfrak{Q}
\circ_{s}^{0\mathbf{Pth}_{\boldsymbol{\mathcal{A}}}}
\mathfrak{P}^{1,\bb{\mathfrak{P}}-1}
\right)\right)\right)
\\=
f_{s}\left(
\mathrm{ip}^{(1,X)@}_{s}\left(
\mathrm{CH}^{(1)}_{s}\left(
\mathfrak{Q}
\right)\right)\right)
\circ_{s}^{0\mathbf{B}}
f_{s}\left(
\mathrm{ip}^{(1,X)@}_{s}\left(
\mathrm{CH}^{(1)}_{s}\left(
\mathfrak{P}^{1,\bb{\mathfrak{P}}-1}
\right)\right)\right);
\end{multline*}
the sixth equality follows from the fact that, since $\mathbf{B}$ a partial many-sorted $\Sigma^{\boldsymbol{\mathcal{A}}}$ in $\mathbf{PAlg}(\boldsymbol{\mathcal{E}}^{\boldsymbol{\mathcal{A}}})$, the $0$-composition is associative; the seventh equality follows from the fact that, by assumption, $f$ is a $\Sigma^{\boldsymbol{\mathcal{A}}}$-homomorphism; the eighth equality follows from the fact that, by Proposition~\ref{PIpCH},  both $\mathrm{ip}^{(1,X)@}_{s}(
\mathrm{CH}^{(1)}_{s}(\mathfrak{P}^{1,\bb{\mathfrak{P}}-1}
))$ and $\mathrm{ip}^{(1,X)@}_{s}(
\mathrm{CH}^{(1)}_{s}(\mathfrak{P}^{0,0}))$ are paths. Thus, the interpretation of the operation symbol for the $0$-composition, $\circ^{0}_{s}$,  in the $\Sigma^{\boldsymbol{\mathcal{A}}}$-algebra $\mathbf{F}_{\Sigma^{\boldsymbol{\mathcal{A}}}}(\mathbf{Pth}_{\boldsymbol{\mathcal{A}}})$ becomes the corresponding interpretation of the $0$-composition $\circ^{0}_{s}$ in the $\Sigma^{\boldsymbol{\mathcal{A}}}$-algebra $\mathbf{Pth}_{\boldsymbol{\mathcal{A}}}$; the ninth equality follows from the fact that, by Definition~\ref{DIp}, $\mathrm{ip}^{(1,X)@}$ is a $\Sigma^{\boldsymbol{\mathcal{A}}}$-homomorphism;  finally, the last equality recovers the value of the Curry-Howard mapping at $\mathfrak{P}$, as we have discussed above.

The case (1.1.2) of $\mathfrak{P}$ being a path of length strictly greater than one containing an echelon on its first step follows.

The case (1.1) of $i=0$ follows.

For case (1.2), i.e., if $i\neq 0$,  since $\bb{\mathfrak{Q}\circ_{s}^{0\mathbf{Pth}_{\boldsymbol{\mathcal{A}}}}\mathfrak{P}}=\bb{\mathfrak{Q}}+\bb{\mathfrak{P}}$, then either (1.2.1) $i\in\bb{\mathfrak{P}}$ or (1.2.2) $i\in [\bb{\mathfrak{P}},\bb{\mathfrak{Q}\circ_{s}^{0\mathbf{Pth}_{\boldsymbol{\mathcal{A}}}}\mathfrak{P}}-1]$.

If (1.2.1), i.e., if $i\neq 0$ and $i\in\bb{\mathfrak{P}}$, then $\mathfrak{P}$ is a path of length strictly greater than one containing an echelon in a step different from the initial one. Then the value of the Curry-Howard mapping at $\mathfrak{P}$ is given by
$$
\mathrm{CH}^{(1)}_{s}\left(
\mathfrak{P}
\right)
=
\mathrm{CH}^{(1)}_{s}\left(
\mathfrak{P}^{i,\bb{\mathfrak{P}}-1}
\right)
\circ_{s}^{0\mathbf{T}_{\Sigma^{\boldsymbol{\mathcal{A}}}}(X)}
\mathrm{CH}^{(1)}_{s}\left(
\mathfrak{P}^{0,i-1}
\right).
$$

Moreover, $(\mathfrak{Q}\circ^{0\mathbf{Pth}_{\boldsymbol{\mathcal{A}}}}_{s}\mathfrak{P})^{i,\bb{\mathfrak{Q}\circ^{0\mathbf{Pth}_{\boldsymbol{\mathcal{A}}}}_{s}\mathfrak{P}}-1}=\mathfrak{Q}\circ_{s}^{0\mathbf{Pth}_{\boldsymbol{\mathcal{A}}}}\mathfrak{P}^{i,\bb{\mathfrak{P}}-1}$ Hence, the value of the Curry-Howard mapping at $\mathfrak{Q}\circ^{0\mathbf{Pth}}_{s}\mathfrak{P}$ is given by
$$
\mathrm{CH}^{(1)}_{s}\left(
\mathfrak{Q}\circ^{0\mathbf{Pth}_{\boldsymbol{\mathcal{A}}}}_{s}\mathfrak{P}
\right)
=
\mathrm{CH}^{(1)}_{s}\left(
\mathfrak{Q}\circ^{0\mathbf{Pth}_{\boldsymbol{\mathcal{A}}}}_{s}\mathfrak{P}
^{i,\bb{\mathfrak{P}}-1}
\right)
\circ_{s}^{0\mathbf{T}_{\Sigma^{\boldsymbol{\mathcal{A}}}}(X)}
\mathrm{CH}^{(1)}_{s}\left(
\mathfrak{P}
^{0,i-1}
\right).
$$

Therefore, the following chain of equalities holds
\begin{flushleft}
$
f_{s}\left(
\mathrm{ip}^{(1,X)@}_{s}\left(
\mathrm{CH}^{(1)}_{s}\left(
\mathfrak{Q}
\circ^{0\mathbf{Pth}_{\boldsymbol{\mathcal{A}}}}_{s}
\mathfrak{P}
\right)\right)\right)
$
\allowdisplaybreaks
\begin{align*}
\qquad&=
f_{s}\left(
\mathrm{ip}^{(1,X)@}_{s}\left(
\mathrm{CH}^{(1)}_{s}\left(
\mathfrak{Q}\circ^{0\mathbf{Pth}_{\boldsymbol{\mathcal{A}}}}_{s}\mathfrak{P}
^{i,\bb{\mathfrak{P}}-1}
\right)
\circ_{s}^{0\mathbf{T}_{\Sigma^{\boldsymbol{\mathcal{A}}}}(X)}
\mathrm{CH}^{(1)}_{s}\left(
\mathfrak{P}
^{0,i-1}
\right)\right)\right)
\tag{1}
\\&=
f_{s}\Big(
\mathrm{ip}^{(1,X)@}_{s}\left(
\mathrm{CH}^{(1)}_{s}\left(
\mathfrak{Q}\circ^{0\mathbf{Pth}_{\boldsymbol{\mathcal{A}}}}_{s}\mathfrak{P}
^{i,\bb{\mathfrak{P}}-1}
\right)\right)
\circ_{s}^{0\mathbf{F}_{\Sigma^{\boldsymbol{\mathcal{A}}}}
(\mathbf{Pth}_{\boldsymbol{\mathcal{A}}})}
\\&\qquad\qquad\qquad\qquad\qquad\qquad\qquad\qquad\qquad
\mathrm{ip}^{(1,X)@}_{s}\left(
\mathrm{CH}^{(1)}_{s}\left(
\mathfrak{P}^{0,i-1}
\right)\right)
\Big)
\tag{2}
\\&=
f_{s}\left(
\mathrm{ip}^{(1,X)@}_{s}\left(
\mathrm{CH}^{(1)}_{s}\left(
\mathfrak{Q}\circ^{0\mathbf{Pth}_{\boldsymbol{\mathcal{A}}}}_{s}\mathfrak{P}
^{i,\bb{\mathfrak{P}}-1}
\right)\right)
\circ_{s}^{0\mathbf{Pth}_{\boldsymbol{\mathcal{A}}}}
\right.
\\&\qquad\qquad\qquad\qquad\qquad\qquad\qquad\qquad\qquad
\left.
\mathrm{ip}^{(1,X)@}_{s}\left(
\mathrm{CH}^{(1)}_{s}\left(
\mathfrak{P}^{0,i-1}
\right)\right)\right)
\tag{3}
\\&=
f_{s}\left(
\mathrm{ip}^{(1,X)@}_{s}\left(
\mathrm{CH}^{(1)}_{s}\left(
\mathfrak{Q}\circ^{0\mathbf{Pth}_{\boldsymbol{\mathcal{A}}}}_{s}\mathfrak{P}
^{i,\bb{\mathfrak{P}}-1}
\right)\right)\right)
\circ_{s}^{0\mathbf{B}}
\\&\qquad\qquad\qquad\qquad\qquad\qquad\qquad\qquad\qquad
f_{s}\left(
\mathrm{ip}^{(1,X)@}_{s}\left(
\mathrm{CH}^{(1)}_{s}\left(
\mathfrak{P}^{0,i-1}
\right)\right)\right)
\tag{4}
\\&=
\left(
f_{s}\left(
\mathrm{ip}^{(1,X)@}_{s}\left(
\mathrm{CH}^{(1)}_{s}\left(
\mathfrak{Q}
\right)\right)\right)
\circ_{s}^{0\mathbf{B}}
f_{s}\left(
\mathrm{ip}^{(1,X)@}_{s}\left(
\mathrm{CH}^{(1)}_{s}\left(
\mathfrak{P}^{i,\bb{\mathfrak{P}}-1}
\right)\right)\right)
\right)
\circ_{s}^{0\mathbf{B}}
\\&\qquad\qquad\qquad\qquad\qquad\qquad\qquad\qquad\qquad
f_{s}\left(
\mathrm{ip}^{(1,X)@}_{s}\left(
\mathrm{CH}^{(1)}_{s}\left(
\mathfrak{P}^{0,i-1}
\right)\right)\right)
\tag{5}
\\&=
f_{s}\left(
\mathrm{ip}^{(1,X)@}_{s}\left(
\mathrm{CH}^{(1)}_{s}\left(
\mathfrak{Q}
\right)\right)\right)
\circ_{s}^{0\mathbf{B}}
\left(
f_{s}\left(
\mathrm{ip}^{(1,X)@}_{s}\left(
\mathrm{CH}^{(1)}_{s}\left(
\mathfrak{P}^{i,\bb{\mathfrak{P}}-1}
\right)\right)\right)
\circ_{s}^{0\mathbf{B}}
\right.
\\&\qquad\qquad\qquad\qquad\qquad\qquad\qquad\qquad
\left.
f_{s}\left(
\mathrm{ip}^{(1,X)@}_{s}\left(
\mathrm{CH}^{(1)}_{s}\left(
\mathfrak{P}^{0,i-1}
\right)\right)\right)
\right)
\tag{6}
\\&=
f_{s}\left(
\mathrm{ip}^{(1,X)@}_{s}\left(
\mathrm{CH}^{(1)}_{s}\left(
\mathfrak{Q}
\right)\right)\right)
\circ_{s}^{0\mathbf{B}}
f_{s}\left(
\mathrm{ip}^{(1,X)@}_{s}\left(
\mathrm{CH}^{(1)}_{s}\left(
\mathfrak{P}^{i,\bb{\mathfrak{P}}-1}
\right)\right)
\circ_{s}^{0\mathbf{Pth}_{\boldsymbol{\mathcal{A}}}}
\right.
\\&\qquad\qquad\qquad\qquad\qquad\qquad\qquad\qquad\qquad
\left.
\mathrm{ip}^{(1,X)@}_{s}\left(
\mathrm{CH}^{(1)}_{s}\left(
\mathfrak{P}^{0,i-1}
\right)\right)\right)
\tag{7}
\\&=
f_{s}\left(
\mathrm{ip}^{(1,X)@}_{s}\left(
\mathrm{CH}^{(1)}_{s}\left(
\mathfrak{Q}
\right)\right)\right)
\circ_{s}^{0\mathbf{B}}
f_{s}\left(
\mathrm{ip}^{(1,X)@}_{s}\left(
\mathrm{CH}^{(1)}_{s}\left(
\mathfrak{P}^{i,\bb{\mathfrak{P}}-1}
\right)\right)
\right.
\\&\qquad\qquad\qquad\qquad\qquad\qquad\quad
\left.
\circ_{s}^{0\mathbf{F}_{\Sigma^{\boldsymbol{\mathcal{A}}}}(
\mathbf{Pth}_{\boldsymbol{\mathcal{A}}})}
\mathrm{ip}^{(1,X)@}_{s}\left(
\mathrm{CH}^{(1)}_{s}\left(
\mathfrak{P}^{0,i-1}
\right)\right)\right)
\tag{8}
\\&=
f_{s}\left(
\mathrm{ip}^{(1,X)@}_{s}\left(
\mathrm{CH}^{(1)}_{s}\left(
\mathfrak{Q}
\right)\right)\right)
\circ_{s}^{0\mathbf{B}}
\\&\qquad\qquad
f_{s}\left(
\mathrm{ip}^{(1,X)@}_{s}\left(
\mathrm{CH}^{(1)}_{s}\left(
\mathfrak{P}^{i,\bb{\mathfrak{P}}-1}
\right)
\circ_{s}^{0\mathbf{T}_{\Sigma^{\boldsymbol{\mathcal{A}}}}(X)}
\mathrm{CH}^{(1)}_{s}\left(
\mathfrak{P}^{0,i-1}
\right)\right)\right)
\tag{9}
\\&=
f_{s}\left(
\mathrm{ip}^{(1,X)@}_{s}\left(
\mathrm{CH}^{(1)}_{s}\left(
\mathfrak{Q}
\right)\right)\right)
\circ_{s}^{0\mathbf{B}}
f_{s}\left(
\mathrm{ip}^{(1,X)@}_{s}\left(
\mathrm{CH}^{(1)}_{s}\left(
\mathfrak{P}
\right)\right)\right).
\tag{10}
\end{align*}
\end{flushleft}

The first equality follows from the previous discussion on the value of the Curry-Howard mapping at $\mathfrak{Q}\circ^{0\mathbf{Pth}_{\boldsymbol{\mathcal{A}}}}_{s}\mathfrak{P}$; the second equality follows from the fact that, by Definition~\ref{DIp}, $\mathrm{ip}^{(1,X)@}$ is a $\Sigma^{\boldsymbol{\mathcal{A}}}$-homomorphism; the third equality follows from the fact that, by Proposition~\ref{PIpCH},  both $\mathrm{ip}^{(1,X)@}_{s}(
\mathrm{CH}^{(1)}_{s}(\mathfrak{Q}
\circ_{s}^{0\mathbf{Pth}_{\boldsymbol{\mathcal{A}}}}\mathfrak{P}^{i,\bb{\mathfrak{P}}-1}
))$ and $\mathrm{ip}^{(1,X)@}_{s}(
\mathrm{CH}^{(1)}_{s}(\mathfrak{P}^{0,i-1}))$ are  paths. Thus, the interpretation of the operation symbol for the $0$-composition, $\circ^{0}_{s}$,  in the $\Sigma^{\boldsymbol{\mathcal{A}}}$-algebra $\mathbf{F}_{\Sigma^{\boldsymbol{\mathcal{A}}}}(\mathbf{Pth}_{\boldsymbol{\mathcal{A}}})$ becomes the corresponding interpretation of the $0$-composition $\circ^{0}_{s}$ in the $\Sigma^{\boldsymbol{\mathcal{A}}}$-algebra $\mathbf{Pth}_{\boldsymbol{\mathcal{A}}}$; the fourth equality follows from the fact that, by assumption, $f$ is a $\Sigma^{\boldsymbol{\mathcal{A}}}$-homomorphism; the fifth equality follows by induction. Let us note that the pair $(\mathfrak{Q}\circ^{0\mathbf{Pth}_{\boldsymbol{\mathcal{A}}}}_{s}\mathfrak{P}^{i,\bb{\mathfrak{P}}-1},s)$ $\prec_{\mathbf{Pth}_{\boldsymbol{\mathcal{A}}}}$-precedes $(\mathfrak{Q}\circ^{0\mathbf{Pth}_{\boldsymbol{\mathcal{A}}}}_{s}\mathfrak{P},s)$. Hence we have that 
\begin{multline*}
f_{s}\left(
\mathrm{ip}^{(1,X)@}_{s}\left(
\mathrm{CH}^{(1)}_{s}\left(
\mathfrak{Q}
\circ_{s}^{0\mathbf{Pth}_{\boldsymbol{\mathcal{A}}}}
\mathfrak{P}^{i,\bb{\mathfrak{P}}-1}
\right)\right)\right)
\\=
f_{s}\left(
\mathrm{ip}^{(1,X)@}_{s}\left(
\mathrm{CH}^{(1)}_{s}\left(
\mathfrak{Q}
\right)\right)\right)
\circ_{s}^{0\mathbf{B}}
f_{s}\left(
\mathrm{ip}^{(1,X)@}_{s}\left(
\mathrm{CH}^{(1)}_{s}\left(
\mathfrak{P}^{i,\bb{\mathfrak{P}}-1}
\right)\right)\right);
\end{multline*}
the sixth equality follows from the fact that, since $\mathbf{B}$ a partial many-sorted $\Sigma^{\boldsymbol{\mathcal{A}}}$ in $\mathbf{PAlg}(\boldsymbol{\mathcal{E}}^{\boldsymbol{\mathcal{A}}})$, the $0$-composition is associative; the seventh equality follows from the fact that, by assumption, $f$ is a $\Sigma^{\boldsymbol{\mathcal{A}}}$-homomorphism; the eighth equality follows from the fact that, by Proposition~\ref{PIpCH},  both $\mathrm{ip}^{(1,X)@}_{s}(
\mathrm{CH}^{(1)}_{s}(\mathfrak{P}^{i,\bb{\mathfrak{P}}-1}
))$ and $\mathrm{ip}^{(1,X)@}_{s}(
\mathrm{CH}^{(1)}_{s}(\mathfrak{P}^{0,i-1}))$ are  paths. Thus, the interpretation of the operation symbol for the $0$-composition, $\circ^{0}_{s}$,  in the $\Sigma^{\boldsymbol{\mathcal{A}}}$-algebra $\mathbf{F}_{\Sigma^{\boldsymbol{\mathcal{A}}}}(\mathbf{Pth}_{\boldsymbol{\mathcal{A}}})$ becomes the corresponding interpretation of the $0$-composition $\circ^{0}_{s}$ in the $\Sigma^{\boldsymbol{\mathcal{A}}}$-algebra $\mathbf{Pth}_{\boldsymbol{\mathcal{A}}}$; the ninth equality follows from the fact that, by Definition~\ref{DIp}, $\mathrm{ip}^{(1,X)@}$ is a $\Sigma^{\boldsymbol{\mathcal{A}}}$-homomorphism;  finally, the last equality recovers the value of the Curry-Howard mapping at $\mathfrak{P}$, as we have discussed above.

The case (1.2.1) $i\in\bb{\mathfrak{P}}$ follows.

If~(1.2.2), i.e., $i\neq 0$ and $i\in[\bb{\mathfrak{P}},\bb{\mathfrak{Q}\circ_{s}^{0\mathbf{Pth}_{\boldsymbol{\mathcal{A}}}}\mathfrak{P}}-1]$ then $\mathfrak{Q}$ is a non $(1,0)$-identity path containing an echelon, whilst $\mathfrak{P}$ is an echelonless path.

We will distinguish three cases according to whether (1.2.2.1) $\mathfrak{Q}$ is an echelon; (1.2.2.2) $\mathfrak{Q}$ is a path of length strictly greater than one containing an echelon on its first step or (1.2.2.3) $\mathfrak{Q}$ is a path of length strictly greater than one containing an echelon on a step different from zero. These cases can be proved following a similar argument to those three cases presented above. We leave the details for the interested reader.

If~(2), i.e., if $\mathfrak{Q}\circ^{0\mathbf{
Pth}_{\boldsymbol{\mathcal{A}}}}_{s}\mathfrak{P}$ is an echelonless path then, regarding the paths $\mathfrak{Q}$ and $\mathfrak{P}$, we have that
\begin{enumerate}
\item[(i)] $\mathfrak{P}$ is an echelonless path;
\item[(ii)] $\mathfrak{Q}$ is an echelonless path.
\end{enumerate}

Since (i), we have by Lemma~\ref{LPthHeadCt}, that there exists a  unique word $\mathbf{s}\in S^{\star}-\{\lambda\}$ and a unique operation symbol $\sigma\in \Sigma_{\mathbf{s},s}$ associated to $\mathfrak{P}$. Let $(\mathfrak{P}_{j})_{j\in\bb{\mathbf{s}}}$ be the family of paths we can extract from $\mathfrak{P}$ in virtue of Lemma~\ref{LPthExtract}. Then, according to Definition~\ref{DCH}, we have that the value of the Curry-Howard mapping at $\mathfrak{P}$ is given by
$$
\mathrm{CH}^{(1)}_{s}\left(
\mathfrak{P}
\right)
=
\sigma^{\mathbf{T}_{\Sigma^{\boldsymbol{\mathcal{A}}}}(X)}\left(\left(
\mathrm{CH}^{(1)}_{s_{j}}\left(
\mathfrak{P}_{j}
\right)
\right)_{j\in\bb{\mathbf{s}}}\right).
$$

Since (ii), we have by Lemma~\ref{LPthHeadCt}, that for the  unique word $\mathbf{s}\in S^{\star}-\{\lambda\}$ and the unique operation symbol $\sigma\in \Sigma_{\mathbf{s},s}$, $\sigma$ is the operation symbol associated to $\mathfrak{Q}$. Note that the operation symbol $\sigma$ is the same as in case (i), since $\mathfrak{Q}\circ^{0\mathbf{Pth}_{\boldsymbol{\mathcal{A}}}}\mathfrak{P}$ is an echelonless path by hypothesis and, thus, head-constant by Lemma~\ref{LPthHeadCt}.

Let $(\mathfrak{Q}_{j})_{j\in\bb{\mathbf{s}}}$ be the family of paths we can extract from $\mathfrak{Q}$ in virtue of Lemma~\ref{LPthExtract}. Then, according to Definition~\ref{DCH}, we have that the value of the Curry-Howard mapping at $\mathfrak{Q}$ is given by
$$
\mathrm{CH}^{(1)}_{s}\left(
\mathfrak{Q}
\right)
=
\sigma^{\mathbf{T}_{\Sigma^{\boldsymbol{\mathcal{A}}}}(X)}\left(\left(
\mathrm{CH}^{(1)}_{s_{j}}\left(
\mathfrak{Q}_{j}
\right)
\right)_{j\in\bb{\mathbf{s}}}\right).
$$

Let $((\mathfrak{Q}\circ^{0\mathbf{Pth}_{\boldsymbol{\mathcal{A}}}}_{s}\mathfrak{P})_{j})_{j\in\bb{\mathbf{s}}}$ be the family of paths in $\mathrm{Pth}_{\boldsymbol{\mathcal{A}},\mathbf{s}}$ which, in virtue of Lemma~\ref{LPthExtract}, we can extract from $\mathfrak{Q}\circ^{0\mathbf{Pth}_{\boldsymbol{\mathcal{A}}}}_{s}\mathfrak{P}$. Then, for every $j\in\bb{\mathbf{s}}$, we have
$$
\left(\mathfrak{Q}\circ^{0\mathbf{Pth}_{\boldsymbol{\mathcal{A}}}}_{s}\mathfrak{P}\right)_{j}
=
\mathfrak{Q}_{j}
\circ^{0\mathbf{Pth}_{\boldsymbol{\mathcal{A}}}}_{s_{j}}
\mathfrak{P}_{j}.
$$

In this case, the value of the Curry-Howard mapping at $\mathfrak{Q}\circ^{0\mathbf{Pth}_{\boldsymbol{\mathcal{A}}}}_{s}\mathfrak{P}$ is given by
$$
\mathrm{CH}^{(1)}_{s}\left(
\mathfrak{Q}\circ^{0\mathbf{Pth}_{\boldsymbol{\mathcal{A}}}}_{s}\mathfrak{P}
\right)=\sigma^{\mathbf{T}_{\Sigma^{\boldsymbol{\mathcal{A}}}}(X)}
\left(\left(
\mathrm{CH}^{(1)}_{s_{j}}\left(
\mathfrak{Q}_{j}\circ^{0\mathbf{Pth}_{\boldsymbol{\mathcal{A}}}}_{s_{j}}\mathfrak{P}_{j}
\right)\right)_{j\in\bb{\mathbf{s}}}\right).
$$

The following chain of equalities holds
\begin{flushleft}
$
f_{s}\left(
\mathrm{ip}^{(1,X)@}_{s}\left(
\mathrm{CH}^{(1)}_{s}\left(
\mathfrak{Q}
\circ^{0\mathbf{Pth}_{\boldsymbol{\mathcal{A}}}}_{s}
\mathfrak{P}
\right)\right)\right)
$
\allowdisplaybreaks
\begin{align*}
&=
f_{s}\left(
\mathrm{ip}^{(1,X)@}_{s}\left(
\sigma^{\mathbf{T}_{\Sigma^{\boldsymbol{\mathcal{A}}}}(X)}
\left(\left(
\mathrm{CH}^{(1)}_{s_{j}}\left(
\mathfrak{Q}_{j}
\circ^{0\mathbf{Pth}_{\boldsymbol{\mathcal{A}}}}_{s_{j}}
\mathfrak{P}_{j}
\right)\right)_{j\in\bb{\mathbf{s}}}
\right)\right)\right)
\tag{1}
\\&=
f_{s}\left(
\sigma^{\mathbf{F}_{\Sigma^{\boldsymbol{\mathcal{A}}}}(
\mathbf{Pth}_{\boldsymbol{\mathcal{A}}})}
\left(\left(
\mathrm{ip}^{(1,X)@}_{s_{j}}\left(
\mathrm{CH}^{(1)}_{s_{j}}\left(
\mathfrak{Q}_{j}
\circ^{0\mathbf{Pth}_{\boldsymbol{\mathcal{A}}}}_{s_{j}}
\mathfrak{P}_{j}
\right)\right)\right)_{j\in\bb{\mathbf{s}}}
\right)\right)
\tag{2}
\\&=
f_{s}\left(
\sigma^{\mathbf{Pth}_{\boldsymbol{\mathcal{A}}}}
\left(\left(
\mathrm{ip}^{(1,X)@}_{s_{j}}\left(
\mathrm{CH}^{(1)}_{s_{j}}\left(
\mathfrak{Q}_{j}
\circ^{0\mathbf{Pth}_{\boldsymbol{\mathcal{A}}}}_{s_{j}}
\mathfrak{P}_{j}
\right)\right)\right)_{j\in\bb{\mathbf{s}}}
\right)\right)
\tag{3}
\\&=
\sigma^{\mathbf{B}}\left(\left(
f_{s_{j}}\left(
\mathrm{ip}^{(1,X)@}_{s_{j}}\left(
\mathrm{CH}^{(1)}_{s_{j}}
\left(
\mathfrak{Q}_{j}
\circ^{0\mathbf{Pth}_{\boldsymbol{\mathcal{A}}}}_{s_{j}}
\mathfrak{P}_{j}
\right)\right)\right)\right)_{j\in\bb{\mathbf{s}}}
\right)
\tag{4}
\\&=
\sigma^{\mathbf{B}}
\left(\left(
f_{s_{j}}\left(
\mathrm{ip}^{(1,X)@}_{s_{j}}\left(
\mathrm{CH}^{(1)}_{s_{j}}\left(
\mathfrak{Q}_{j}
\right)\right)\right)
\circ^{0\mathbf{B}}_{s_{j}}
\right.\right.
\\&\qquad\qquad\qquad\qquad\qquad\qquad \qquad\qquad
\left.\left.
f_{s_{j}}\left(
\mathrm{ip}^{(1,X)@}_{s_{j}}\left(
\mathrm{CH}^{(1)}_{s_{j}}\left(
\mathfrak{P}_{j}
\right)\right)\right)
\right)_{j\in\bb{\mathbf{s}}}\right)
\tag{5}
\\&=
\sigma^{\mathbf{B}}
\left(\left(
f_{s_{j}}\left(
\mathrm{ip}^{(1,X)@}_{s_{j}}\left(
\mathrm{CH}^{(1)}_{s_{j}}\left(
\mathfrak{Q}_{j}
\right)\right)\right)\right)_{j\in\bb{\mathbf{s}}}\right)
\circ^{0\mathbf{B}}_{s}
\\&\qquad\qquad\qquad\qquad\qquad\qquad\qquad
\sigma^{\mathbf{B}}
\left(\left(
f_{s_{j}}\left(
\mathrm{ip}^{(1,X)@}_{s_{j}}\left(
\mathrm{CH}^{(1)}_{s_{j}}\left(
\mathfrak{P}_{j}
\right)\right)\right)\right)_{j\in\bb{\mathbf{s}}}
\right)
\tag{6}
\\&=
f_{s}\left(
\sigma^{\mathbf{Pth}_{\boldsymbol{\mathcal{A}}}}
\left(\left(
\mathrm{ip}^{(1,X)@}_{s_{j}}\left(
\mathrm{CH}^{(1)}_{s_{j}}\left(
\mathfrak{Q}_{j}
\right)\right)\right)_{j\in\bb{\mathbf{s}}}
\right)\right)
\circ^{0\mathbf{B}}_{s}
\\&\qquad\qquad\qquad\qquad\qquad\qquad\quad
f_{s}\left(
\sigma^{\mathbf{Pth}_{\boldsymbol{\mathcal{A}}}}
\left(\left(
\mathrm{ip}^{(1,X)@}_{s_{j}}\left(
\mathrm{CH}^{(1)}_{s_{j}}\left(
\mathfrak{P}_{j}
\right)\right)\right)_{j\in\bb{\mathbf{s}}}
\right)\right)
\tag{7}
\\&=
f_{s}\left(
\sigma^{\mathbf{F}_{\Sigma^{\boldsymbol{\mathcal{A}}}}(
\mathbf{Pth}_{\boldsymbol{\mathcal{A}}})}
\left(\left(
\mathrm{ip}^{(1,X)@}_{s_{j}}\left(
\mathrm{CH}^{(1)}_{s_{j}}\left(
\mathfrak{Q}_{j}
\right)\right)\right)_{j\in\bb{\mathbf{s}}}
\right)\right)
\circ^{0\mathbf{B}}_{s}
\\&\qquad\qquad\qquad\qquad\qquad
f_{s}\left(
\sigma^{\mathbf{F}_{\Sigma^{\boldsymbol{\mathcal{A}}}}(
\mathbf{Pth}_{\boldsymbol{\mathcal{A}}})}
\left(\left(
\mathrm{ip}^{(1,X)@}_{s_{j}}\left(
\mathrm{CH}^{(1)}_{s_{j}}\left(
\mathfrak{P}_{j}
\right)\right)\right)_{j\in\bb{\mathbf{s}}}
\right)\right)
\tag{8}
\\
&=
f_{s}\left(
\mathrm{ip}^{(1,X)@}_{s}\left(
\sigma^{\mathbf{T}_{\Sigma^{\boldsymbol{\mathcal{A}}}}(X)}
\left(\left(
\mathrm{CH}^{(1)}_{s_{j}}\left(
\mathfrak{Q}_{j}
\right)
\right)_{j\in\bb{\mathbf{s}}}
\right)\right)\right)
\circ^{0\mathbf{B}}_{s}
\\&\qquad\qquad\qquad\qquad\qquad\qquad
f_{s}\left(
\mathrm{ip}^{(1,X)@}_{s}\left(
\sigma^{\mathbf{T}_{\Sigma^{\boldsymbol{\mathcal{A}}}}(X)}
\left(\left(
\mathrm{CH}^{(1)}_{s_{j}}\left(
\mathfrak{P}_{j}
\right)\right)_{j\in\bb{\mathbf{s}}}
\right)\right)\right)
\tag{9}
\\&=
f_{s}\left(
\mathrm{ip}^{(1,X)@}_{s}\left(
\mathrm{CH}^{(1)}_{s}\left(
\mathfrak{Q}
\right)\right)\right)
\circ^{0\mathbf{B}}_{s}
f_{s}\left(
\mathrm{ip}^{(1,X)@}_{s}\left(
\mathrm{CH}^{(1)}_{s}\left(
\mathfrak{P}
\right)\right)\right).
\tag{10}
\end{align*}
\end{flushleft}

The first equality holds because it consists, simply, in unravelling the definition of the Curry-Howard mapping at $\mathfrak{Q}\circ^{0\mathbf{Pth}_{\boldsymbol{\mathcal{A}}}}_{s}\mathfrak{P}$ as we have discussed above; with regard to the second equality, we apply that $\mathrm{ip}^{(1,X)@}$ is a $\Sigma^{\boldsymbol{\mathcal{A}}}$-homomorphism, according to Definition~\ref{DIp}; the third equation follows from the fact that, in virtue of Proposition~\ref{PIpCH},  for every $j\in\bb{\mathbf{s}}$, the element
$\mathrm{ip}^{(1,X)@}_{s_{j}}(
\mathrm{CH}^{(1)}_{s_{j}}(\mathfrak{Q}_{j}\circ_{s_{j}}^{0\mathbf{Pth}_{\boldsymbol{\mathcal{A}}}}\mathfrak{P}_{j}
))$  is a path. Thus, the interpretation of the operation symbol $\sigma$,  in the $\Sigma^{\boldsymbol{\mathcal{A}}}$-algebra $\mathbf{F}_{\Sigma^{\boldsymbol{\mathcal{A}}}}(\mathbf{Pth}_{\boldsymbol{\mathcal{A}}})$ becomes the corresponding interpretation of $\sigma$ in the $\Sigma^{\boldsymbol{\mathcal{A}}}$-algebra $\mathbf{Pth}_{\boldsymbol{\mathcal{A}}}$; the fourth equation holds because, by assumption, $f$ is  a $\Sigma^{\boldsymbol{\mathcal{A}}}$-homomorphism; the fifth equation holds by induction, since, for every $j\in\bb{\mathbf{s}}$, the pair $(\mathfrak{Q}_{j}\circ^{0\mathbf{Pth}_{\boldsymbol{\mathcal{A}}}}_{s_{j}}\mathfrak{P}_{j},s_{j})$ $\prec_{\mathbf{Pth}_{\boldsymbol{\mathcal{A}}}}$-precedes $(\mathfrak{Q}\circ^{0\mathbf{Pth}_{\boldsymbol{\mathcal{A}}}}_{s}\mathfrak{P},s)$. Therefore the statement holds for these paths extracted in virtue of Lemma~\ref{LPthExtract}, that is
\begin{multline*}
f_{s_{j}}\left(
\mathrm{ip}^{(1,X)@}_{s_{j}}\left(
\mathrm{CH}^{(1)}_{s_{j}}\left(
\mathfrak{Q}_{j}
\circ^{0\mathbf{Pth}_{\boldsymbol{\mathcal{A}}}}_{s_{j}}
\mathfrak{P}_{j}
\right)\right)\right)
\\=
f_{s_{j}}\left(
\mathrm{ip}^{(1,X)@}_{s_{j}}\left(
\mathrm{CH}^{(1)}_{s_{j}}\left(
\mathfrak{Q}_{j}
\right)\right)\right)
\circ^{0\mathbf{B}}_{s_{j}}
f_{s_{j}}\left(
\mathrm{ip}^{(1,X)@}_{s_{j}}\left(
\mathrm{CH}^{(1)}_{s_{j}}\left(
\mathfrak{P}_{j}
\right)\right)\right);
\end{multline*}
the sixth equality holds because the many-sorted partial $\Sigma^{\boldsymbol{\mathcal{A}}}$-algebra $\mathbf{B}$ is in $\mathbf{PAlg}(\boldsymbol{\mathcal{E}}^{\boldsymbol{\mathcal{A}}})$, thus, by Axiom~\ref{DVarA8}, the $0$-composition  $\circ^{0}_{s}$ is compatible with the realization of the operation symbol $\sigma$ in $\mathbf{B}$; the seventh equality holds because, by assumption, $f$ is a $\Sigma^{\boldsymbol{\mathcal{A}}}$-homomorphism; the eighth equality follows from the fact that, in virtue of Proposition~\ref{PIpCH},  for every $j\in\bb{\mathbf{s}}$, 
both $\mathrm{ip}^{(1,X)@}_{s_{j}}(
\mathrm{CH}^{(1)}_{s_{j}}(\mathfrak{Q}_{j}
))$ and $\mathrm{ip}^{(1,X)@}_{s_{j}}(
\mathrm{CH}^{(1)}_{s_{j}}(\mathfrak{P}_{j}
))$  are  paths. Thus, the interpretation of the operation symbol $\sigma$,  in the $\Sigma^{\boldsymbol{\mathcal{A}}}$-algebra $\mathbf{F}_{\Sigma^{\boldsymbol{\mathcal{A}}}}(\mathbf{Pth}_{\boldsymbol{\mathcal{A}}})$ becomes the corresponding interpretation of $\sigma$ in the $\Sigma^{\boldsymbol{\mathcal{A}}}$-algebra $\mathbf{Pth}_{\boldsymbol{\mathcal{A}}}$; the ninth equality holds because $\mathrm{ip}^{(1,X)@}$ is a $\Sigma^{\boldsymbol{\mathcal{A}}}$-homomorphism, according to Definition~\ref{DIp}; finally, the last equality simply unravels the definition of the Curry-Howard mapping at the paths $\mathfrak{Q}$ and $\mathfrak{P}$, respectively.

The case (2) follows.

It follows that Equation~\ref{PVarKerE2} holds, i.e., the mapping $\mathrm{pr}^{
\equiv^{[1]}
}\circ\mathrm{ip}^{(1,X)@}\circ\mathrm{CH}^{(1)}$ is compatible with the $0$-composition.

All in all, we  conclude that $\mathrm{pr}^{
\equiv^{[1]}
}\circ\mathrm{ip}^{(1,X)@}\circ\mathrm{CH}^{(1)}$ is a $\Sigma^{\boldsymbol{\mathcal{A}}}$-homomorphism.

It remains to prove that 
\[
\mathrm{Ker}(\mathrm{CH}^{(1)})
\subseteq \mathrm{Ker}(\mathrm{pr}^{
\equiv^{[1]}
}\circ\mathrm{ip}^{(1,X)@}\circ\mathrm{CH}^{(1)}).\]

Let $s$ be a sort in $S$ and let $\mathfrak{Q}$, $\mathfrak{P}$ be two paths in $\mathrm{Pth}_{\boldsymbol{\mathcal{A}},s}$ satisfying that $(\mathfrak{Q},\mathfrak{P})\in\mathrm{Ker}(\mathrm{CH})_{s}$.

The following chain of equalities holds
\[
\mathrm{pr}^{
\equiv^{[1]}
}_{s}\left(\mathrm{ip}^{(1,X)@}_{s}\left(\mathrm{CH}^{(1)}_{s}\left(
\mathfrak{Q}
\right)\right)\right)
\\=
\mathrm{pr}^{
\equiv^{[1]}
}_{s}\left(\mathrm{ip}^{(1,X)@}_{s}\left(\mathrm{CH}^{(1)}_{s}\left(
\mathfrak{P}
\right)\right)\right).
\]

It follows directly from the fact that $\mathrm{CH}^{(1)}_{s}(\mathfrak{Q})=\mathrm{CH}^{(1)}_{s}(\mathfrak{P})$. Consequently, \[
\mathrm{Ker}(\mathrm{CH}^{(1)})
\subseteq \mathrm{Ker}(\mathrm{pr}^{
\equiv^{[1]}
}\circ\mathrm{ip}^{(1,X)@}\circ\mathrm{CH}^{(1)}).\]

This completes the proof of Proposition~\ref{PVarKer}.
\end{proof}

\begin{restatable}{remark}{RQPUniv}
\label{RQPUniv} From Proposition~\ref{PVarKer} and taking into account the universal property of the quotient, we can assert that there exists a unique $\Sigma^{\boldsymbol{\mathcal{A}}}$-homomorphism, that we will denote by $(\mathrm{pr}^{\equiv^{[1]}}
\circ
\mathrm{ip}^{(1,X)@}
\circ
\mathrm{CH}^{(1)})^{\natural}$, i.e.,
\[
(\mathrm{pr}^{\equiv^{[1]}}
\circ
\mathrm{ip}^{(1,X)@}
\circ
\mathrm{CH}^{(1)})^{\natural}
\colon 
\left[
\mathbf{Pth}_{\boldsymbol{\mathcal{A}}}
\right]
\mor
\mathbf{T}_{\boldsymbol{\mathcal{E}}^{\boldsymbol{\mathcal{A}}}}\left(
\mathbf{Pth}_{\boldsymbol{\mathcal{A}}}
\right)
\]
such that 
$
(\mathrm{pr}^{\equiv^{[1]}}
\circ
\mathrm{ip}^{(1,X)@}
\circ
\mathrm{CH}^{(1)})^{\natural}
\circ 
\mathrm{pr}^{\mathrm{Ker}(\mathrm{CH}^{(1)})}
=
\mathrm{pr}^{\equiv^{[1]}}
\circ
\mathrm{ip}^{(1,X)@}
\circ
\mathrm{CH}^{(1)}.
$
\end{restatable}

\begin{restatable}{proposition}{PQEta}
\label{PQEta}
The  diagram of Figure~\ref{FQEta} commutes, i.e., the following equality holds
$$
\mathrm{pr}^{\equiv^{[1]}}
\circ
\mathrm{ip}^{(1,X)@}
\circ
\mathrm{CH}^{(1)}
=
\eta^{([1],\mathbf{Pth}_{\boldsymbol{\mathcal{A}}})}.
$$
\end{restatable}

\begin{figure}
\begin{tikzpicture}
[ACliment/.style={-{To [angle'=45, length=5.75pt, width=4pt, round]}
}, scale=0.8]
\node[] (P3) 		at 	(0,-4) 	[] 	{$\mathbf{Pth}_{\boldsymbol{\mathcal{A}}}$};
\node[] (TXQ) 	at 	(6,-4) 	[] 	{$\mathbf{PT}_{\boldsymbol{\mathcal{A}}}
(X)$};
\node[] (SchP) 		at 	(6,-6) 	[] 	{$\mathbf{Sch}_{\boldsymbol{\mathcal{E}}^{\boldsymbol{\mathcal{A}}}}(\mathbf{Pth}_{\boldsymbol{\mathcal{A}}})$};
\node[] (TEP2) 	at 	(6,-8) 	[] 	{$\mathbf{T}_{\boldsymbol{\mathcal{E}}^{\boldsymbol{\mathcal{A}}}}
(\mathbf{Pth}_{\boldsymbol{\mathcal{A}}})$};

\draw[ACliment]  (TXQ) 	to node [right]	
{$\mathrm{ip}^{(1,X)@}$} (SchP);
\draw[ACliment]  (SchP) 	to node [right]	
{$\mathrm{pr}^{\equiv^{[1]}}$}  (TEP2);
\draw[ACliment, bend right=10]  (P3) 	to node [below left]	
{$\eta^{([1],\mathbf{Pth}_{\boldsymbol{\mathcal{A}}})}$} (TEP2);

\draw[ACliment]  (P3) 	to node [above]	
{$\mathrm{CH}^{(1)}$} (TXQ);
\end{tikzpicture}
\caption{The maps in Proposition~\ref{PQEta}.}
\label{FQEta}
\end{figure}

\begin{proof}
Let $s$ be a sort in $S$ and $\mathfrak{P}$ a path in $\mathrm{Pth}_{\boldsymbol{\mathcal{A}},s}$. We need to prove that
$$
\left[\mathfrak{P}\right]_{\equiv^{[1]}_{s}}=
\left[
\mathrm{ip}^{(1,X)@}_{s}\left(
\mathrm{CH}^{(1)}_{s}\left(
\mathfrak{P}
\right)\right)\right]_{\equiv^{[1]}_{s}}.
$$
But this is equivalent to prove that, for every many-sorted partial $\Sigma^{\boldsymbol{\mathcal{A}}}$-algebra $\mathbf{B}$ in $\mathbf{PAlg}(\boldsymbol{\mathcal{E}}^{\boldsymbol{\mathcal{A}}})$ and every $\Sigma^{\boldsymbol{\mathcal{A}}}$-homomorphism $f\colon\mathbf{Pth}_{\boldsymbol{\mathcal{A}}}\mor\mathbf{B}$ we have
$$
f^{\mathrm{Sch}}_{s}\left(
\mathfrak{P}
\right)=
f^{\mathrm{Sch}}_{s}\left(
\mathrm{ip}^{(1,X)@}_{s}\left(
\mathrm{CH}^{(1)}_{s}\left(
\mathfrak{P}
\right)\right)\right).
$$
However, by Proposition~\ref{PIpCH}, the element $\mathrm{ip}^{(1,X)@}_{s}(\mathrm{CH}^{(1)}_{s}(\mathfrak{P}))$ is a path. Therefore, the last equation reduces to the equation
$$
f_{s}\left(
\mathfrak{P}\right)
=f_{s}\left(
\mathrm{ip}^{(1,X)@}_{s}\left(
\mathrm{CH}^{(1)}_{s}\left(
\mathfrak{P}
\right)\right)\right).
$$

We prove this statement by Artinian induction on $(\coprod\mathrm{Pth}_{\boldsymbol{\mathcal{A}}},\leq_{\mathbf{Pth}_{\boldsymbol{\mathcal{A}}}})$.

\textsf{Base step of the Artinian induction}.

Let $(\mathfrak{P},s)$ be a minimal element in $(\coprod\mathrm{Pth}_{\boldsymbol{\mathcal{A}}}, \leq_{\mathbf{Pth}_{\boldsymbol{\mathcal{A}}}})$. Then, by Proposition~\ref{PMinimal}, the path $\mathfrak{P}$ is either a $(1,0)$-identity path on simple term, or an echelon. In any case, by Corollary~\ref{CIpCHOneStep}, if $\mathfrak{P}$ is a $(1,0)$-identity path or a one-step path, then the following equation holds
$$
\mathfrak{P}=\mathrm{ip}^{(1,X)@}_{s}\left(
\mathrm{CH}^{(1)}_{s}\left(
\mathfrak{P}\right)\right).
$$
Therefore, the statement holds for paths of length $0$ or $1$.

\textsf{Inductive step of the Artinian induction}.

Let $(\mathfrak{P},s)$ be a non-minimal element in $(\coprod\mathrm{Pth}_{\boldsymbol{\mathcal{A}}}, \leq_{\mathbf{Pth}_{\boldsymbol{\mathcal{A}}}})$. Let us suppose that, for every sort $t\in S$ and every path $\mathfrak{Q}\in\mathrm{Pth}_{\boldsymbol{\mathcal{A}},t}$, if $(\mathfrak{Q},t)$ $<_{\mathbf{Pth}_{\boldsymbol{\mathcal{A}}}}$-precedes $(\mathfrak{P},s)$, then the statement holds for $\mathfrak{Q}$, i.e., we have that
$$
f_{t}\left(
\mathfrak{Q}
\right)=f_{t}\left(
\mathrm{ip}^{(1,X)@}_{t}\left(
\mathrm{CH}^{(1)}_{t}\left(
\mathfrak{Q}
\right)\right)\right).
$$

Let us recall that, by the discussion presented in the base case, we can also assume that $\mathfrak{P}$ is not a $(1,0)$-identity path. Since $(\mathfrak{P},s)$ is a non-minimal element in $(\coprod\mathrm{Pth}_{\boldsymbol{\mathcal{A}}}, \leq_{\mathbf{Pth}_{\boldsymbol{\mathcal{A}}}})$, we have, by Lemma~\ref{LOrdI}, that $\mathfrak{P}$ is either~(1) a path of length strictly greater than one containing at least one echelon or~(2) an echelonless path. 

If~(1), then let $i\in \bb{\mathfrak{P}}$ be the first index for which the one-step subpath $\mathfrak{P}^{i,i}$ of $\mathfrak{P}$ is an echelon. We distinguish the cases (1.1) $i=0$ and (1.2) $i>0$.

If~(1.1), i.e., if $\mathfrak{P}$ is a path of length strictly greater than one having an echelon on its first step,  then the value of the Curry-Howard mapping at $\mathfrak{P}$ is given by
$$
\mathrm{CH}^{(1)}_{s}\left(
\mathfrak{P}
\right)=
\mathrm{CH}^{(1)}_{s}\left(
\mathfrak{P}^{1,\bb{\mathfrak{P}}-1}
\right)
\circ_{s}^{0\mathbf{T}_{\Sigma^{\boldsymbol{\mathcal{A}}}}(X)}
\mathrm{CH}^{(1)}_{s}\left(
\mathfrak{P}^{0,0}
\right).
$$

Thus, the following chain of equalities holds
\begin{flushleft}
$f_{s}\left(
\mathrm{ip}^{(1,X)@}_{s}\left(
\mathrm{CH}^{(1)}_{s}\left(
\mathfrak{P}
\right)\right)\right)$
\allowdisplaybreaks
\begin{align*}
&=
f_{s}\left(
\mathrm{ip}^{(1,X)@}_{s}\left(
\mathrm{CH}^{(1)}_{s}\left(
\mathfrak{P}^{1,\bb{\mathfrak{P}}-1}
\right)
\circ_{s}^{0\mathbf{T}_{\Sigma^{\boldsymbol{\mathcal{A}}}}(X)}
\mathrm{CH}^{(1)}_{s}\left(
\mathfrak{P}^{0,0}
\right)
\right)\right)
\tag{1}
\\&=
f_{s}\left(
\mathrm{ip}^{(1,X)@}_{s}\left(
\mathrm{CH}^{(1)}_{s}\left(
\mathfrak{P}^{1,\bb{\mathfrak{P}}-1}
\right)\right)
\circ_{s}^{0\mathbf{F}_{\Sigma^{\boldsymbol{\mathcal{A}}}}
(\mathbf{Pth}_{\boldsymbol{\mathcal{A}}})}
\mathrm{ip}^{(1,X)@}_{s}\left(
\mathrm{CH}^{(1)}_{s}\left(
\mathfrak{P}^{0,0}
\right)\right)\right)
\tag{2}
\\&=
f_{s}\left(
\mathrm{ip}^{(1,X)@}_{s}\left(
\mathrm{CH}^{(1)}_{s}\left(
\mathfrak{P}^{1,\bb{\mathfrak{P}}-1}
\right)\right)
\circ_{s}^{0\mathbf{Pth}_{\boldsymbol{\mathcal{A}}}}
\mathrm{ip}^{(1,X)@}_{s}\left(
\mathrm{CH}^{(1)}_{s}\left(
\mathfrak{P}^{0,0}
\right)\right)\right)
\tag{3}
\\&=
f_{s}\left(
\mathrm{ip}^{(1,X)@}_{s}\left(
\mathrm{CH}^{(1)}_{s}\left(
\mathfrak{P}^{1,\bb{\mathfrak{P}}-1}
\right)\right)\right)
\circ_{s}^{0\mathbf{B}}
f_{s}\left(
\mathrm{ip}^{(1,X)@}_{s}\left(
\mathrm{CH}^{(1)}_{s}\left(
\mathfrak{P}^{0,0}
\right)\right)\right)
\tag{4}
\\&=
f_{s}\left(
\mathfrak{P}^{1,\bb{\mathfrak{P}}-1}
\right)
\circ_{s}^{0\mathbf{B}}
f_{s}\left(
\mathfrak{P}^{0,0}
\right)
\tag{5}
\\&=
f_{s}\left(
\mathfrak{P}^{1,\bb{\mathfrak{P}}-1}
\circ_{s}^{0\mathbf{\mathbf{Pth}_{\boldsymbol{\mathcal{A}}}}}
\mathfrak{P}^{0,0}
\right)
\tag{6}
\\&=
f_{s}\left(
\mathfrak{P}
\right).
\tag{7}
\end{align*}
\end{flushleft}

In the just stated chain of equalities, the first equality holds because it consists, simply, in unravelling the definition of the Curry-Howard mapping at $\mathfrak{P}$; the second equality holds because $\mathrm{ip}^{(1,X)@}$ is a $\Sigma^{\boldsymbol{\mathcal{A}}}$-homomorphism, according to Definition~\ref{DIp}; the third equation follows from the fact that, in virtue of Proposition~\ref{PIpCH},  both $\mathrm{ip}^{(1,X)@}_{s}(
\mathrm{CH}^{(1)}_{s}(\mathfrak{P}^{1,\bb{\mathfrak{P}}-1}))$ and $\mathrm{ip}^{(1,X)@}_{s}(
\mathrm{CH}^{(1)}_{s}(\mathfrak{P}^{0,0}))$ are paths. Thus, the interpretation of the operation symbol for the $0$-composition, $\circ^{0}_{s}$,  in the $\Sigma^{\boldsymbol{\mathcal{A}}}$-algebra $\mathbf{F}_{\Sigma^{\boldsymbol{\mathcal{A}}}}(\mathbf{Pth}_{\boldsymbol{\mathcal{A}}})$ becomes the corresponding interpretation of the $0$-composition $\circ^{0}_{s}$ in the $\Sigma^{\boldsymbol{\mathcal{A}}}$-algebra $\mathbf{Pth}_{\boldsymbol{\mathcal{A}}}$; the fourth equality holds because $f$ is a $\Sigma^{\boldsymbol{\mathcal{A}}}$-homomorphism;
the fifth equality holds by induction, since the pairs $(\mathfrak{P}^{1,\bb{\mathfrak{P}}-1},s)$ and $(\mathfrak{P}^{0,0},s)$ $\prec_{\mathbf{Pth}_{\boldsymbol{\mathcal{A}}}}$-precede $(\mathfrak{P},s)$; the fifth equation holds since $f$ is a $\Sigma^{\boldsymbol{\mathcal{A}}}$-homomorphism; and the last equality holds because in it we, simply, are recovering the original path $\mathfrak{P}$ as the composition of its $0$-constituents.

If~(1.2), i.e., if $\mathfrak{P}$ is a path of length strictly greater than one having its first echelon on a step $i\in\bb{\mathfrak{P}}$ different from the initial one, then the value of the Curry-Howard mapping at $\mathfrak{P}$ is given by
$$
\mathrm{CH}^{(1)}_{s}\left(
\mathfrak{P}
\right)=
\mathrm{CH}^{(1)}_{s}\left(
\mathfrak{P}^{i,\bb{\mathfrak{P}}-1}
\right)
\circ_{s}^{0\mathbf{T}_{\Sigma^{\boldsymbol{\mathcal{A}}}}(X)}
\mathrm{CH}^{(1)}_{s}\left(
\mathfrak{P}^{0,i-1}
\right).
$$

Thus, the following chain of equalities holds
\begin{flushleft}
$f_{s}\left(
\mathrm{ip}^{(1,X)@}_{s}\left(
\mathrm{CH}^{(1)}_{s}\left(
\mathfrak{P}
\right)\right)\right)$
\allowdisplaybreaks
\begin{align*}
&=
f_{s}\left(
\mathrm{ip}^{(1,X)@}_{s}\left(
\mathrm{CH}^{(1)}_{s}\left(
\mathfrak{P}^{i,\bb{\mathfrak{P}}-1}
\right)
\circ_{s}^{0\mathbf{T}_{\Sigma^{\boldsymbol{\mathcal{A}}}}(X)}
\mathrm{CH}^{(1)}_{s}\left(
\mathfrak{P}^{0,i-1}
\right)\right)\right)
\tag{1}
\\&=
f_{s}\left(
\mathrm{ip}^{(1,X)@}_{s}\left(
\mathrm{CH}^{(1)}_{s}\left(
\mathfrak{P}^{i,\bb{\mathfrak{P}}-1}
\right)\right)
\circ_{s}^{0\mathbf{F}_{\Sigma^{\boldsymbol{\mathcal{A}}}}
(\mathbf{Pth}_{\boldsymbol{\mathcal{A}}})}
\mathrm{ip}^{(1,X)@}_{s}\left(
\mathrm{CH}^{(1)}_{s}\left(
\mathfrak{P}^{0,i-1}
\right)\right)\right)
\tag{2}
\\&=
f_{s}\left(
\mathrm{ip}^{(1,X)@}_{s}\left(
\mathrm{CH}^{(1)}_{s}\left(
\mathfrak{P}^{i,\bb{\mathfrak{P}}-1}
\right)\right)
\circ_{s}^{0\mathbf{Pth}_{\boldsymbol{\mathcal{A}}}}
\mathrm{ip}^{(1,X)@}_{s}\left(
\mathrm{CH}^{(1)}_{s}\left(
\mathfrak{P}^{0,i-1}
\right)\right)\right)
\tag{3}
\\&=
f_{s}\left(
\mathrm{ip}^{(1,X)@}_{s}\left(
\mathrm{CH}^{(1)}_{s}\left(
\mathfrak{P}^{i,\bb{\mathfrak{P}}-1}
\right)\right)\right)
\circ_{s}^{0\mathbf{B}}
f_{s}\left(
\mathrm{ip}^{(1,X)@}_{s}\left(
\mathrm{CH}^{(1)}_{s}\left(
\mathfrak{P}^{0,i-1}
\right)\right)\right)
\tag{4}
\\&=
f_{s}\left(
\mathfrak{P}^{i,\bb{\mathfrak{P}}-1}
\right)
\circ_{s}^{0\mathbf{B}}
f_{s}\left(
\mathfrak{P}^{0,i-1}
\right)
\tag{5}
\\&=
f_{s}\left(
\mathfrak{P}^{i,\bb{\mathfrak{P}}-1}
\circ_{s}^{0\mathbf{\mathbf{Pth}_{\boldsymbol{\mathcal{A}}}}}
\mathfrak{P}^{0,i-1}
\right)
\tag{6}
\\&=
f_{s}\left(
\mathfrak{P}
\right).
\tag{7}
\end{align*}
\end{flushleft}

In the just stated chain of equalities, the first equality holds because it consists, simply, in unravelling the definition of the Curry-Howard mapping at $\mathfrak{P}$; the second equality holds because $\mathrm{ip}^{(1,X)@}$ is a $\Sigma^{\boldsymbol{\mathcal{A}}}$-homomorphism, according to Definition~\ref{DIp}; the third equation follows from the fact that, in virtue of Proposition~\ref{PIpCH},  both $\mathrm{ip}^{(1,X)@}_{s}(
\mathrm{CH}^{(1)}_{s}(\mathfrak{P}^{i,\bb{\mathfrak{P}}-1}))$ and $\mathrm{ip}^{(1,X)@}_{s}(
\mathrm{CH}^{(1)}_{s}(\mathfrak{P}^{0,i-1}))$ are proper paths. Thus, the interpretation of the operation symbol for the $0$-composition, $\circ^{0}_{s}$,  in the $\Sigma^{\boldsymbol{\mathcal{A}}}$-algebra $\mathbf{F}_{\Sigma^{\boldsymbol{\mathcal{A}}}}(\mathbf{Pth}_{\boldsymbol{\mathcal{A}}})$ becomes the corresponding interpretation of the $0$-composition $\circ^{0}_{s}$ in the $\Sigma^{\boldsymbol{\mathcal{A}}}$-algebra $\mathbf{Pth}_{\boldsymbol{\mathcal{A}}}$; the fourth equality holds because $f$ is a $\Sigma^{\boldsymbol{\mathcal{A}}}$-homomorphism;
the fifth equality holds by induction, since the pairs $(\mathfrak{P}^{i,\bb{\mathfrak{P}}-1},s)$ and $(\mathfrak{P}^{0,i-1},s)$ $\prec_{\mathbf{Pth}_{\boldsymbol{\mathcal{A}}}}$-precede $(\mathfrak{P},s)$; the fifth equation holds since $f$ is a $\Sigma^{\boldsymbol{\mathcal{A}}}$-homomorphism; finally, the last equality holds because in it we, simply, are recovering the original path $\mathfrak{P}$ as the composition of its $0$-constituents.

If~(2), i.e., if $\mathfrak{P}$ is an echelonless path then, according to Lemma~\ref{LPthHeadCt}, there exists a unique word $\mathbf{s}\in S^{\star}-\{\lambda\}$ and a unique operation symbol $\sigma\in \Sigma_{\mathbf{s},s}$ associated to $\mathfrak{P}$. Let $(\mathfrak{P}_{j})_{j\in\bb{\mathbf{s}}}$ be the family of paths in $\mathrm{Pth}_{\boldsymbol{\mathcal{A}},\mathbf{s}}$ which, in virtue of Lemma~\ref{LPthExtract}, we can extract from $\mathfrak{P}$.  In this case, the value of the Curry-Howard mapping at $\mathfrak{P}$ is given by
$$
\mathrm{CH}^{(1)}_{s}\left(\mathfrak{P}\right)=
\sigma^{\mathbf{T}_{\Sigma^{\boldsymbol{\mathcal{A}}}}(X)}
\left(\left(
\mathrm{CH}^{(1)}_{s_{j}}\left(
\mathfrak{P}_{j}
\right)\right)_{j\in\bb{\mathbf{s}}}
\right).
$$

Thus, the following chain of equalities holds
\begin{flushleft}
$f_{s}\left(
\mathrm{ip}^{(1,X)@}_{s}\left(
\mathrm{CH}^{(1)}_{s}\left(
\mathfrak{P}
\right)\right)\right)$
\allowdisplaybreaks
\begin{align*}
\qquad
&=
f_{s}\left(
\mathrm{ip}^{(1,X)@}_{s}\left(
\sigma^{\mathbf{T}_{\Sigma^{\boldsymbol{\mathcal{A}}}}(X)}
\left(\left(
\mathrm{CH}^{(1)}_{s_{j}}\left(
\mathfrak{P}_{j}
\right)\right)_{j\in\bb{\mathbf{s}}}
\right)\right)\right)
\tag{1}
\\&=
f_{s}\left(
\sigma^{\mathbf{F}_{\Sigma^{\boldsymbol{\mathcal{A}}}}(
\mathbf{Pth}_{\boldsymbol{\mathcal{A}}})}
\left(\left(
\mathrm{ip}^{(1,X)@}_{s_{j}}\left(
\mathrm{CH}^{(1)}_{s_{j}}\left(
\mathfrak{P}_{j}
\right)\right)\right)_{j\in\bb{\mathbf{s}}}
\right)\right)
\tag{2}
\\&=
f_{s}\left(
\sigma^{\mathbf{Pth}_{\boldsymbol{\mathcal{A}}}}
\left(\left(
\mathrm{ip}^{(1,X)@}_{s_{j}}\left(
\mathrm{CH}^{(1)}_{s_{j}}\left(
\mathfrak{P}_{j}
\right)\right)\right)_{j\in\bb{\mathbf{s}}}
\right)\right)
\tag{3}
\\&=
\sigma^{\mathbf{B}}
\left(\left(
f_{s_{j}}\left(
\mathrm{ip}^{(1,X)@}_{s_{j}}\left(
\mathrm{CH}^{(1)}_{s_{j}}\left(
\mathfrak{P}_{j}
\right)\right)\right)\right)_{j\in\bb{\mathbf{s}}}
\right)
\tag{4}
\\&=
\sigma^{\mathbf{B}}
\left(\left(
f_{s_{j}}\left(
\mathfrak{P}_{j}
\right)\right)_{j\in\bb{\mathbf{s}}}
\right)
\tag{5}
\\&=
f_{s}\left(
\sigma^{\mathbf{Pth}_{\boldsymbol{\mathcal{A}}}}
\left(\left(
\mathfrak{P}_{j}
\right)_{j\in\bb{\mathbf{s}}}
\right)\right).
\tag{6}
\end{align*}
\end{flushleft}

In the just stated chain of equalities, the first equality holds because it consists, simply, in unravelling the definition of the Curry-Howard mapping at $\mathfrak{P}$; the second equality holds because $\mathrm{ip}^{(1,X)@}$ is a $\Sigma^{\boldsymbol{\mathcal{A}}}$-homomorphism, according to Definition~\ref{DIp}; the third equation follows from the fact that, in virtue of Proposition~\ref{PIpCH},  for every $j\in\bb{\mathbf{s}}$, the element
$\mathrm{ip}^{(1,X)@}_{s_{j}}(
\mathrm{CH}^{(1)}_{s_{j}}(\mathfrak{P}_{j}
))$  is a  path. Thus, the interpretation of the operation symbol $\sigma$,  in the $\Sigma^{\boldsymbol{\mathcal{A}}}$-algebra $\mathbf{F}_{\Sigma^{\boldsymbol{\mathcal{A}}}}(\mathbf{Pth}_{\boldsymbol{\mathcal{A}}})$ becomes the corresponding interpretation of $\sigma$ in the $\Sigma^{\boldsymbol{\mathcal{A}}}$-algebra $\mathbf{Pth}_{\boldsymbol{\mathcal{A}}}$; the fourth equality holds because $f$ is a $\Sigma^{\boldsymbol{\mathcal{A}}}$-homomorphism; the fifth equality holds by induction, since, for every $j\in\bb{\mathbf{s}}$, the pair $((\mathfrak{P}_{j},s_{j}),(\mathfrak{P},s))$ is in $\prec_{\coprod{\mathrm{Pth}_{\boldsymbol{\mathcal{A}}}}}$; finally, the last equality holds because $f$  is a $\Sigma^{\boldsymbol{\mathcal{A}}}$-homomorphism.

It remains to prove that, for an echelonless path $\mathfrak{P}$, the following equality holds
$$
f_{s}\left(
\mathfrak{P}
\right)=f_{s}\left(
\sigma^{\mathbf{Pth}_{\boldsymbol{\mathcal{A}}}}
\left(\left(
\mathfrak{P}_{j}
\right)_{j\in\bb{\mathbf{s}}}
\right)\right).
$$

Let $\mathbf{c}$ be the unique word of $S^{\star}$ for which $\mathfrak{P}$ is a path in $\mathrm{Pth}_{\mathbf{c},\boldsymbol{\mathcal{A}},s}$. Since $\mathfrak{P}$ is an echelonless path, we have, by Lemma~\ref{LPthExtract}, that there exists a partition of $\mathbf{c}$ into a family $(\mathbf{c}_{j})_{j\in\bb{\mathbf{s}}}$ in $(S^{\star})^{\bb{\mathbf{s}}}$ such that, for every $j\in\bb{\mathbf{s}}$, $\mathfrak{P}_{j}$ is a path in $\mathrm{Pth}_{\mathbf{c}_{j},\boldsymbol{\mathcal{A}},s_{j}}$.
Moreover, by Proposition~\ref{PPthRecons}, we have that $\mathfrak{P}$ can be univocally represented as a composition of one-step paths as
\begin{align*}
\mathfrak{P}&=
\mathfrak{P}^{\bb{\mathbf{c}}-1,\bb{\mathbf{c}}-1}
\circ^{0\mathbf{Pth}_{\boldsymbol{\mathcal{A}}}}_{s}
\cdots
\circ^{0\mathbf{Pth}_{\boldsymbol{\mathcal{A}}}}_{s}
\mathfrak{P}^{0,0}.
\tag{E1}
\label{PQEtaE1}
\end{align*}

Moreover, since $\mathfrak{P}$ is an echelonless path, we have, by Lemma~\ref{LPthExtract}, that there exists a partition of $\mathbf{c}$ into a family $(\mathbf{c}_{j})_{j\in\bb{\mathbf{s}}}$ in $(S^{\star})^{\bb{\mathbf{s}}}$ such that, for every $j\in\bb{\mathbf{s}}$, $\mathfrak{P}_{j}$ is a  path in $\mathrm{Pth}_{\mathbf{c}_{j},\boldsymbol{\mathcal{A}},s_{j}}$.

Let us note that, for every $i\in \bb{\mathbf{c}}$, the one-step  subpath $\mathfrak{P}^{i,i}$ of $\mathfrak{P}$ is an echelonless  path associated to the operation symbol $\sigma$. Now, recovering the terminology used in Lemma~\ref{LPthExtract}, the index $i$ must have type $(j_{i},k_{i})$, for some $j_{i}\in\bb{\mathbf{s}}$ and some $k_{i}\in\bb{\mathbf{c}_{j_{i}}}$, meaning that the unique  translation appearing on the one-step  path $\mathfrak{P}^{i,i}$ has its derived translation occurring at position $j_{i}$, this fact being the $k_{i}$-th time that occurs, i.e.,  there exists a family of terms $(P_{i,l})_{l\in j_{i}}\in \prod_{l\in j_{i}}\mathrm{T}_{\Sigma}(X)_{s_{l}}$ and a family of terms $(P_{i,l})_{l\in \bb{\mathbf{s}}-(j_{i}+1)}\in \prod_{l\in \bb{\mathbf{s}}-(j_{i}+1)}\mathrm{T}_{\Sigma}(X)_{s_{l}}$ and a  translation $T'_{i}\in \mathrm{Tl}_{c_{i}}(\mathbf{T}_{\Sigma}(X))_{s_{j_{i}}}$ such that
\[
T_{i}=\sigma^{\mathbf{T}_{\Sigma}(X)}\left(
P_{i,0},
\cdots,
P_{i,j_{i}-1},
T'_{i},
P_{i,j_{i}+1},
\cdots,
P_{i,\bb{\mathbf{s}}-1}
\right).
\]

Let $((\mathfrak{P}^{i,i})_{j})_{j\in\bb{\mathbf{s}}}$ be the family of  paths that we can extract from $\mathfrak{P}^{i,i}$ in virtue of Lemma~\ref{LPthExtract}. What we obtain, following the proof of Lemma~\ref{LPthExtract}, is the following family of  paths
\allowdisplaybreaks
\begin{align*}
\left(\mathfrak{P}^{i,i}\right)_{0}&=
\mathrm{ip}^{(1,0)\sharp}_{s_{0}}\left(
P_{i,0}
\right);
\\
\vdots&\qquad\vdots
\\
\left(\mathfrak{P}^{i,i}\right)_{j_{i}-1}&=
\mathrm{ip}^{(1,0)\sharp}_{s_{j_{i}-1}}\left(
P_{i,j_{i}-1}
\right);
\\
\left(\mathfrak{P}^{i,i}\right)_{j_{i}}&=
\left(\mathfrak{P}_{j_{i}}\right)^{k_{i},k_{i}};
\\
\left(\mathfrak{P}^{i,i}\right)_{j_{i}+1}&=
\mathrm{ip}^{(1,0)\sharp}_{s_{j_{i}+1}}\left(
P_{i,j_{i}+1}
\right);
\\
\vdots&\qquad\vdots
\\
\left(\mathfrak{P}^{i,i}\right)_{\bb{\mathbf{s}}-1}&=
\mathrm{ip}^{(1,0)\sharp}_{s_{\bb{\mathbf{s}}-1}}\left(
P_{i,\bb{\mathbf{s}}-1}
\right).
\end{align*}

That is, we recover $(1,0)$-identity  paths on every index in $\bb{\mathbf{s}}$ different from $j_{i}$, whilst on index $j_{i}$ we recover, precisely, the $k_{i}$-th step of the $j_{i}$-th path that we can extract from $\mathfrak{P}$, i.e., $(\mathfrak{P}_{j_{i}})^{k_{i},k_{i}}$.

Two questions regarding the families of paths $((\mathfrak{P}^{i,i}_{j})_{j\in \bb{\mathbf{s}}})_{i\in \bb{\mathbf{c}}}$ arise. On the first hand, let us note that, since $\mathfrak{P}^{i,i}$ is a one-step echelonless  path, it follows, by Corollary~\ref{CUStep}, that
\begin{align*}
\mathfrak{P}^{i,i}
&=
\sigma^{\mathbf{Pth}_{\boldsymbol{\mathcal{A}}}}
\left(
\left(
\left(\mathfrak{P}^{i,i}\right)_{j}
\right)_{j\in\bb{\mathbf{s}}}
\right).
\tag{E2}
\label{PQEtaE2}
\end{align*}

Moreover, following Proposition~\ref{PPthRecons}, we can reconstruct, for every $j\in\bb{\mathbf{s}}$, the path $\mathfrak{P}_{j}$, as the $0$-composition of its one-step subpaths as follows
$$
\mathfrak{P}_{j}
=
\left(\mathfrak{P}_{j}\right)^{\bb{\mathbf{c}_{j}}-1, \bb{\mathbf{c}_{j}}-1}
\circ^{0\mathbf{Pth}_{\boldsymbol{\mathcal{A}}}}_{s_{j}}
\cdots
\circ^{0\mathbf{Pth}_{\boldsymbol{\mathcal{A}}}}_{s_{j}}
\left(\mathfrak{P}_{j}\right)^{0,0}.
$$

The following equation is also a valid description of $\mathfrak{P}_{j}$ as a  $0$-composition
\begin{align*}
\mathfrak{P}_{j}&=
\left(\mathfrak{P}^{\bb{\mathbf{c}}-1,\bb{\mathbf{c}}-1}\right)_{j}
\circ^{0\mathbf{Pth}_{\boldsymbol{\mathcal{A}}}}_{s_{j}}
\cdots
\circ^{0\mathbf{Pth}_{\boldsymbol{\mathcal{A}}}}_{s_{j}}
\left(\mathfrak{P}^{0,0}\right)_{j}.
\tag{E3}
\label{PQEtaE3}
\end{align*}

In this last equation we are $0$-composing all the  paths at position $j$ that we can extract from the one-step subpaths $\mathfrak{P}^{i,i}$, for every $i\in\bb{\mathbf{c}}$. Note that, for those indexes $i\in\bb{\mathbf{c}}$, that are not of type $j$, the  path $(\mathfrak{P}^{i,i})_{j}$ will simply be a $(1,0)$-identity path that acts as a neutral element for the $0$-composition.

All in all, the following chain of equalities holds
\allowdisplaybreaks
\begin{align*}
f_{s}\left(
\mathfrak{P}
\right)&=
f_{s}\left(
\mathfrak{P}^{\bb{\mathbf{c}}-1,\bb{\mathbf{c}}-1}
\circ^{0\mathbf{Pth}_{\boldsymbol{\mathcal{A}}}}_{s}
\cdots
\circ^{0\mathbf{Pth}_{\boldsymbol{\mathcal{A}}}}_{s}
\mathfrak{P}^{0,0}
\right)
\tag{1}
\\&=
f_{s}\left(
\mathfrak{P}^{\bb{\mathbf{c}}-1,\bb{\mathbf{c}}-1}
\right)
\circ^{1\mathbf{B}}_{s}
\cdots
\circ^{1\mathbf{B}}_{s}
f_{s}
\Big(
\mathfrak{P}^{0,0}
\Big)
\tag{2}
\\&=
f_{s}\left(
\sigma^{\mathbf{Pth}_{\boldsymbol{\mathcal{A}}}}
\left(
\left(
\left(
\mathfrak{P}^{\bb{\mathbf{c}}-1,\bb{\mathbf{c}}-1}
\right)_{j}
\right)_{j\in\bb{\mathbf{s}}}
\right)\right)
\circ^{0\mathbf{B}}_{s}
\cdots
\circ^{0\mathbf{B}}_{s}
\\&
\qquad\qquad\qquad\qquad\qquad\qquad\qquad\qquad
f_{s}\left(
\sigma^{\mathbf{Pth}_{\boldsymbol{\mathcal{A}}}}
\left(
\left(
\left(
\mathfrak{P}^{0,0}
\right)_{j}
\right)_{j\in\bb{\mathbf{s}}}
\right)\right)
\tag{3}
\\&=
\sigma^{\mathbf{B}}
\left(
\left(
f_{s_{j}}\left(
\left(
\mathfrak{P}^{\bb{\mathbf{c}}-1,\bb{\mathbf{c}}-1}
\right)_{j}
\right)
\right)_{j\in\bb{\mathbf{s}}}\right)
\circ^{0\mathbf{B}}_{s}
\cdots
\circ^{0\mathbf{B}}_{s}
\\&
\qquad\qquad\qquad\qquad\qquad\qquad\qquad\qquad\qquad
\sigma^{\mathbf{B}}
\left(
\left(
f_{s_{j}}\left(
\left(
\mathfrak{P}^{0,0}
\right)_{j}
\right)
\right)_{j\in\bb{\mathbf{s}}}\right)
\tag{4}
\\&=
\sigma^{\mathbf{B}}
\left(\left(
f_{s_{j}}\left(
\left(
\mathfrak{P}^{\bb{\mathbf{c}}-1,\bb{\mathbf{c}}-1}
\right)_{j}
\right)
\circ^{0\mathbf{B}}_{s_{j}}
\cdots 
\circ^{0\mathbf{B}}_{s_{j}}
f_{s_{j}}\left(
\left(
\mathfrak{P}^{0,0}
\right)_{j}
\right)
\right)_{j\in\bb{\mathbf{s}}}
\right)
\tag{5}
\\
&=
\sigma^{\mathbf{B}}
\left(\left(
f_{s_{j}}\left(
\left(\mathfrak{P}^{\bb{\mathbf{c}}-1,\bb{\mathbf{c}}-1}\right)_{j}
\circ^{0\mathbf{Pth}_{\boldsymbol{\mathcal{A}}}}_{s_{j}}
\cdots
\circ^{0\mathbf{Pth}_{\boldsymbol{\mathcal{A}}}}_{s_{j}}
\left(\mathfrak{P}^{0,0}\right)_{j}
\right)
\right)_{j\in\bb{\mathbf{s}}}
\right)
\tag{6}
\\
&=
\sigma^{\mathbf{B}}
\left(\left(
f_{s_{j}}\left(
\mathfrak{P}_{j}
\right)\right)_{j\in\bb{\mathbf{s}}}
\right)
\tag{7}
\\&=
f_{s}\left(\sigma^{\mathbf{Pth}_{\boldsymbol{\mathcal{A}}}}
\left(\left(
\mathfrak{P}_{j}
\right)_{j\in\bb{\mathbf{s}}}
\right)\right).
\tag{8}
\end{align*}

In the just stated chain of equalities, the first equality holds because $\mathfrak{P}$ can be decomposed into one-step  paths, as shown in Equation~\ref{PQEtaE1}; the second equality holds because $f$ is a $\Sigma^{\boldsymbol{\mathcal{A}}}$-homomorphism; the third equation follows
from Equation~\ref{PQEtaE2}; the fourth equality holds because $f$ is a $\Sigma^{\boldsymbol{\mathcal{A}}}$-homomorphism; the fifth equality holds because the many-sorted partial $\Sigma^{\boldsymbol{\mathcal{A}}}$-algebra $\mathbf{B}$ is in $\mathbf{PAlg}(\boldsymbol{\mathcal{E}}^{\boldsymbol{\mathcal{A}}})$ thus, by Axiom~\ref{DVarA8}, the $0$-composition,  $\circ^{0}_{s}$, is compatible with the interpretation of the operation $\sigma$ in $\mathbf{B}$; the sixth equality holds because $f$ is a $\Sigma^{\boldsymbol{\mathcal{A}}}$-homomorphism; the seventh equality holds according to Equation~\ref{PQEtaE3}; finally, the last equality holds because $f$ is a $\Sigma^{\boldsymbol{\mathcal{A}}}$-homomorphism.

This completes Case~(2).

This completes the proof.
\end{proof}

With the above propositions we can introduce the main result of this chapter.

\begin{restatable}{theorem}{TPthFree}
\label{TPthFree} 
The many-sorted partial $\Sigma^{\boldsymbol{\mathcal{A}}}$-algebras $[\mathbf{Pth}_{\boldsymbol{\mathcal{A}}}]$ and $\mathbf{T}_{\boldsymbol{\mathcal{E}}^{\boldsymbol{\mathcal{A}}}}(\mathbf{Pth}_{\boldsymbol{\mathcal{A}}})$ are isomorphic.
\end{restatable}
\begin{proof}
To follow the proof, the reader is urged to consider the diagram given in Figure~\ref{FBigPth}. By Corollary~\ref{CVarPr}, the $\Sigma^{\boldsymbol{\mathcal{A}}}$-homomorphism $\mathrm{pr}^{\mathrm{Ker}(\mathrm{CH}^{(1)})}$ from $\mathbf{Pth}_{\boldsymbol{\mathcal{A}}}$ to $[\mathbf{Pth}_{\boldsymbol{\mathcal{A}}}]$ extends to a $\Sigma^{\boldsymbol{\mathcal{A}}}$-homomorphism, $\mathrm{pr}^{\mathrm{Ker}(\mathrm{CH}^{(1)})\mathsf{p}}$ from $\mathbf{T}_{\boldsymbol{\mathcal{E}}^{\boldsymbol{\mathcal{A}}}}(\mathbf{Pth}_{\boldsymbol{\mathcal{A}}})$ to $[\mathbf{Pth}_{\boldsymbol{\mathcal{A}}}]$.  By Remark~\ref{RQPUniv}, the mapping $
(\mathrm{pr}^{\equiv^{[1]}}
\circ
\mathrm{ip}^{(1,X)@}
\circ
\mathrm{CH}^{(1)})^{\natural}
$ from $
[
\mathbf{Pth}_{\boldsymbol{\mathcal{A}}}
]
$ to $
\mathbf{T}_{\boldsymbol{\mathcal{E}}^{\boldsymbol{\mathcal{A}}}}
(
\mathbf{Pth}_{\boldsymbol{\mathcal{A}}}
)
$
is a $\Sigma^{\boldsymbol{\mathcal{A}}}$-homomorphism. 

Thus, the composition
$$
(\mathrm{pr}^{\equiv^{[1]}}
\circ
\mathrm{ip}^{(1,X)@}
\circ
\mathrm{CH}^{(1)})^{\natural}
\circ
\mathrm{pr}^{\mathrm{Ker}(\mathrm{CH}^{(1)})\mathsf{p}}
$$
is a $\Sigma^{\boldsymbol{\mathcal{A}}}$-endomorphism of $\mathbf{T}_{\boldsymbol{\mathcal{E}}^{\boldsymbol{\mathcal{A}}}}(\mathbf{Pth}_{\boldsymbol{\mathcal{A}}})$.

Moreover, this composition satisfies the following chain of equalities
\begin{flushleft}
$(\mathrm{pr}^{\equiv^{[1]}}
\circ
\mathrm{ip}^{(1,X)@}
\circ
\mathrm{CH}^{(1)})^{\natural}
\circ
\mathrm{pr}^{\mathrm{Ker}(\mathrm{CH}^{(1)})\mathsf{p}}
\circ
\eta^{([1],\mathbf{Pth}_{\boldsymbol{\mathcal{A}}})}$
\allowdisplaybreaks
\begin{align*}
\qquad
&=
(\mathrm{pr}^{\equiv^{[1]}}
\circ
\mathrm{ip}^{(1,X)@}
\circ
\mathrm{CH}^{(1)})^{\natural}
\circ
\mathrm{pr}^{\mathrm{Ker}(\mathrm{CH}^{(1)})}
\tag{1}
\\
&=
\mathrm{pr}^{\equiv^{[1]}}
\circ
\mathrm{ip}^{(1,X)@}
\circ
\mathrm{CH}^{(1)}
\tag{2}
\\
&=
\eta^{([1],\mathbf{Pth}_{\boldsymbol{\mathcal{A}}})}.
\tag{3}
\end{align*}
\end{flushleft}

In the just stated chain of equalities, the first equality holds by Corollary~\ref{CVarPr}; the second equality follows from Remark~\ref{RQPUniv}; finally the last equality holds by Proposition~\ref{PQEta}.

Let us recall that the identity mapping $\mathrm{id}^{\mathbf{T}_{\boldsymbol{\mathcal{E}}^{\boldsymbol{\mathcal{A}}}}(\mathbf{Pth}_{\boldsymbol{\mathcal{A}}})}$ is also a $\Sigma^{\boldsymbol{\mathcal{A}}}$-endomorphism of $\mathbf{T}_{\boldsymbol{\mathcal{E}}^{\boldsymbol{\mathcal{A}}}}(\mathbf{Pth}_{\boldsymbol{\mathcal{A}}})$ satisfying that
$$
\mathrm{id}^{\mathbf{T}_{\boldsymbol{\mathcal{E}}^{\boldsymbol{\mathcal{A}}}}(\mathbf{Pth}_{\boldsymbol{\mathcal{A}}})}
\circ
\eta^{([1],\mathbf{Pth}_{\boldsymbol{\mathcal{A}}})}
=\eta^{([1],\mathbf{Pth}_{\boldsymbol{\mathcal{A}}})}.
$$

Hence, by the universal property of $\mathbf{T}_{\boldsymbol{\mathcal{E}}^{\boldsymbol{\mathcal{A}}}}(\mathbf{Pth}_{\boldsymbol{\mathcal{A}}})$, we have 
$$
(\mathrm{pr}^{\equiv^{[1]}}
\circ
\mathrm{ip}^{(1,X)@}
\circ
\mathrm{CH}^{(1)})^{\natural}
\circ
\mathrm{pr}^{\mathrm{Ker}(\mathrm{CH}^{(1)})\mathsf{p}}
=
\mathrm{id}^{\mathbf{T}_{\boldsymbol{\mathcal{E}}^{\boldsymbol{\mathcal{A}}}}(\mathbf{Pth}_{\boldsymbol{\mathcal{A}}})}.
$$

Therefore, $(\mathrm{pr}^{\equiv^{[1]}}
\circ
\mathrm{ip}^{(1,X)@}
\circ
\mathrm{CH}^{(1)})^{\natural}$ is a $\Sigma^{\boldsymbol{\mathcal{A}}}$-isomorphism from  $[\mathbf{Pth}_{\boldsymbol{\mathcal{A}}}]$ to $\mathbf{T}_{\boldsymbol{\mathcal{E}}^{\boldsymbol{\mathcal{A}}}}(\mathbf{Pth}_{\boldsymbol{\mathcal{A}}})$ whose inverse is $\mathrm{pr}^{\mathrm{Ker}(\mathrm{CH}^{(1)})\mathsf{p}}$.

This finishes the proof.
\end{proof}

\begin{restatable}{corollary}{CPTFree}
\label{CPTFree} 
$[\mathbf{PT}_{\boldsymbol{\mathcal{A}}}]$ is a partial $\Sigma^{\boldsymbol{\mathcal{A}}}$-algebra in $\mathbf{PAlg}(\boldsymbol{\mathcal{E}}^{\boldsymbol{\mathcal{A}}})$. Indeed, the partial $\Sigma^{\boldsymbol{\mathcal{A}}}$-algebras $[\mathbf{PT}_{\boldsymbol{\mathcal{A}}}]$ and $\mathbf{T}_{\boldsymbol{\mathcal{E}}^{\boldsymbol{\mathcal{A}}}}(\mathbf{Pth}_{\boldsymbol{\mathcal{A}}})$ are isomorphic.
\end{restatable}

\begin{figure}
\begin{center}
\begin{tikzpicture}
[ACliment/.style={-{To [angle'=45, length=5.75pt, width=4pt, round]}
}, scale=0.6]
\node[] (P3) 		at 	(0,-4) 	[] 	{$\mathbf{Pth}_{\boldsymbol{\mathcal{A}}}$};
\node[] (TEP) 	at 	(6,0) 	[] 	{$\mathbf{T}_{\boldsymbol{\mathcal{E}}^{\boldsymbol{\mathcal{A}}}}(\mathbf{Pth}_{\boldsymbol{\mathcal{A}}})$};
\node[] (PK) 		at 	(6,-4) 	[] 	{$[\mathbf{Pth}_{\boldsymbol{\mathcal{A}}}]$};
\node[] (TEP2) 	at 	(6,-8) 	[] 	{$\mathbf{T}_{\boldsymbol{\mathcal{E}}^{\boldsymbol{\mathcal{A}}}}
(\mathbf{Pth}_{\boldsymbol{\mathcal{A}}})$};
\draw[ACliment, bend left]  (P3) 	to node [above left]	
{$\eta^{([1],\mathbf{Pth}_{\boldsymbol{\mathcal{A}}})}$} (TEP);
\draw[ACliment, bend right]  (P3) 	to node [below left]	
{$\eta^{([1],\mathbf{Pth}_{\boldsymbol{\mathcal{A}}})}$} (TEP2);
\draw[ACliment]  (P3) 	to node [above]	
{$\mathrm{pr}^{\mathrm{Ker}(\mathrm{CH}^{(1)})}$} (PK);
\draw[ACliment]  (PK) 	to node [midway, fill=white]	
{$(\mathrm{pr}^{\equiv^{[1]}}
\circ
\mathrm{ip}^{(1,X)@}
\circ
\mathrm{CH}^{(1)})^{\natural}$} (TEP2);
\draw[ACliment]  (TEP) 	to node [right]	
{$\mathrm{pr}^{\mathrm{Ker}(\mathrm{CH}^{(1)})\mathsf{p}}$} (PK);
\draw[ACliment, rounded corners] (TEP.east)
--
 ($(TEP.east)+(1.75,0)$)
to node [right]	
{$\mathrm{id}^{\mathbf{T}_{\boldsymbol{\mathcal{E}}^{\boldsymbol{\mathcal{A}}}}(\mathbf{Pth}_{\boldsymbol{\mathcal{A}}})}$} 
($(TEP2.east)+(1.75,0)$)
-- (TEP2.east);
\end{tikzpicture}
\end{center}
\caption{The maps in Theorem~\ref{TPthFree}.}
\label{FBigPth}
\end{figure}
\chapter{
\texorpdfstring
{First-order translations for path terms}
{First-order translations}
}\label{S1L}

In this chapter we introduce the notions of elementary first-order translation and first-order translation. We show that first-order translations are well-behaved, in the sense that if two path terms $P$ and $Q$ in $\mathrm{PT}_{\boldsymbol{\mathcal{A}}, s}$, for $s\in S$, are such that $\mathrm{ip}^{(1,X)@}_{s}(P)$ and $\mathrm{ip}^{(1,X)@}_{s}(Q)$ have the same $(0,1)$-source and $(0,1)$-target, then the terms that result from applying a first-order translation to $P$ and $Q$ are also path terms and, when they are converted into paths by means of $\mathrm{ip}^{(1,X)@}$, then they still have the same $(0,1)$-source and $(0,1)$-target. Moreover, we show that, for two path terms $P$ and $Q$ in $\mathrm{PT}_{\boldsymbol{\mathcal{A}}, s}$ if they are $\Theta^{[1]}_{s}$-related, then their respective translations are also $\Theta^{[1]}_{s}$-related.


We next define for the many-sorted partial $\Sigma^{\boldsymbol{\mathcal{A}}}$-algebra $\mathbf{T}_{\Sigma^{\boldsymbol{\mathcal{A}}}}(X)$ the concepts of elementary first-order translation and of first-order translation respect to it.

\begin{restatable}{definition}{DUETrans}
\label{DUETrans} 
\index{translation!first-order!elementary}
Let $t$ be a sort in $S$. We will denote by $\mathrm{Etl}_{t}(\mathbf{T}_{\Sigma^{\boldsymbol{\mathcal{A}}}}(X))$ the subset $(\mathrm{Etl}_{t}(\mathbf{T}_{\Sigma^{\boldsymbol{\mathcal{A}}}}(X))_{s})_{s\in S}$ of $(\mathrm{Hom}(\mathrm{T}_{\Sigma^{\boldsymbol{\mathcal{A}}}}(X)_{t},\mathrm{T}_{\Sigma^{\boldsymbol{\mathcal{A}}}}(X)_{s}))_{s\in S}$ defined, for every $s\in S$, as follows: for every mapping $T^{(1)}\in \mathrm{Hom}(\mathrm{T}_{\Sigma^{\boldsymbol{\mathcal{A}}}}(X)_{t},\mathrm{T}_{\Sigma^{\boldsymbol{\mathcal{A}}}}(X)_{s})$, $T^{(1)}\in \mathrm{Etl}_{t}(\mathbf{T}_{\Sigma^{\boldsymbol{\mathcal{A}}}}(X))_{s}$ if and only if one of the following conditions holds
\begin{enumerate}
\item There is a word $\mathbf{s}\in S^{\star}-\{\lambda\}$, an index $k\in\bb{\mathbf{s}}$, an operation symbol $\sigma\in \Sigma_{\mathbf{s},s}$, a family of paths $(\mathfrak{P}_{j})_{j\in k}\in \prod_{j\in k}\mathrm{Pth}_{\boldsymbol{\mathcal{A}},s_{j}}$ and a family of paths $(\mathfrak{P}_{l})_{l\in \bb{\mathbf{s}}-(k+1)}\in \prod_{l\in \bb{\mathbf{s}}-(k+1)}\mathrm{Pth}_{\boldsymbol{\mathcal{A}},s_{l}}$ (recall that $k+1=\{0,\cdots, k\}$ and that $\bb{\mathbf{s}}-(k+1)=\{k+1,\cdots, \bb{\mathbf{s}}-1\}$) such that $s_{k}=t$ and, for every $P\in\mathrm{T}_{\Sigma^{\boldsymbol{\mathcal{A}}}}(X)_{t}$
\begin{multline*}
T^{(1)}(P)=\sigma^{\mathbf{T}_{\Sigma^{\boldsymbol{\mathcal{A}}}}(X)}\left(
\mathrm{CH}^{(1)}_{s_{0}}\left(\mathfrak{P}_{0}\right),
\cdots,
\mathrm{CH}^{(1)}_{s_{k-1}}\left(\mathfrak{P}_{k-1}\right),
\right.
\\
\left.
P,
\mathrm{CH}^{(1)}_{s_{k+1}}\left(\mathfrak{P}_{k+1}\right),
\cdots,
\mathrm{CH}^{(1)}_{s_{\bb{\mathbf{s}}-1}}\left(\mathfrak{P}_{\bb{\mathbf{s}}-1}\right)
\right);
\end{multline*}
\item It holds that $s=t$ and there is a path $\mathfrak{P}\in\mathrm{Pth}_{\boldsymbol{\mathcal{A}},s}$ such that, for every $P\in\mathrm{T}_{\Sigma^{\boldsymbol{\mathcal{A}}}}(X)_{s}$,
\begin{align*}
T^{(1)}(P)&=P\circ^{0\mathbf{T}_{\Sigma^{\boldsymbol{\mathcal{A}}}}(X)}_{s}\mathrm{CH}^{(1)}_{s}\left(\mathfrak{P}\right)
&
\mbox{or}
&&
T^{(1)}(P)&=\mathrm{CH}^{(1)}_{s}\left(\mathfrak{P}\right)\circ^{0\mathbf{T}_{\Sigma^{\boldsymbol{\mathcal{A}}}}(X)}_{s}P.
\end{align*}
\end{enumerate}

We will sometimes use the following presentations of first-order elementary translations, adding an underlined space to denote where the variable will be placed:
\begin{multline*}
T^{(1)}=\sigma^{\mathbf{T}_{\Sigma^{\boldsymbol{\mathcal{A}}}}(X)}\left(
\mathrm{CH}^{(1)}_{s_{0}}\left(\mathfrak{P}_{0}\right),
\cdots,
\mathrm{CH}^{(1)}_{s_{k-1}}\left(\mathfrak{P}_{k-1}\right),
\right.
\\
\left.
\underline{\quad},
\mathrm{CH}^{(1)}_{s_{k+1}}\left(\mathfrak{P}_{k+1}\right),
\cdots,
\mathrm{CH}^{(1)}_{s_{\bb{\mathbf{s}}-1}}\left(\mathfrak{P}_{\bb{\mathbf{s}}-1}\right)
\right);
\tag{1}
\end{multline*}
\begin{align*}
T^{(1)}&=\underline{\quad}\circ^{0\mathbf{T}_{\Sigma^{\boldsymbol{\mathcal{A}}}}(X)}_{s}\mathrm{CH}^{(1)}_{s}\left(\mathfrak{P}\right)
&
\mbox{or}
&&
T^{(1)}&=\mathrm{CH}^{(1)}_{s}\left(\mathfrak{P}\right)\circ^{0\mathbf{T}_{\Sigma^{\boldsymbol{\mathcal{A}}}}(X)}_{s}\underline{\quad}.
\tag{2}
\end{align*}
In the first case we will say that $T^{(1)}$ is a \emph{first-order elementary translation of type $\sigma$}, while in the second case we will say that $T^{(1)}$ is a \emph{first-order elementary translation of type $\circ^{0}_{s}$}. We will call the elements of $\mathrm{Etl}_{t}(\mathbf{T}_{\Sigma^{\boldsymbol{\mathcal{A}}}}(X))_{s}$ the \emph{first-order $t$-elementary translations of sort $s$ for $\mathbf{T}_{\Sigma^{\boldsymbol{\mathcal{A}}}}(X)$}.
\end{restatable}

\begin{restatable}{definition}{DUTrans}
\label{DUTrans} 
\index{translation!first-order}
Let $t$ be a sort in $S$. We will denote by $\mathrm{Tl}_{t}(\mathbf{T}_{\Sigma^{\boldsymbol{\mathcal{A}}}}(X))$ the subset $(\mathrm{Tl}_{t}(\mathbf{T}_{\Sigma^{\boldsymbol{\mathcal{A}}}}(X))_{s})_{s\in S}$ of $(\mathrm{Hom}(\mathrm{T}_{\Sigma^{\boldsymbol{\mathcal{A}}}}(X)_{t},\mathrm{T}_{\Sigma^{\boldsymbol{\mathcal{A}}}}(X)_{s}))_{s\in S}$ defined, for every $s\in S$, as follows: for every mapping $T^{(1)}\in \mathrm{Hom}(\mathrm{T}_{\Sigma^{\boldsymbol{\mathcal{A}}}}(X)_{t},\mathrm{T}_{\Sigma^{\boldsymbol{\mathcal{A}}}}(X)_{s})$, $T^{(1)}\in \mathrm{Tl}_{t}(\mathbf{T}_{\Sigma^{\boldsymbol{\mathcal{A}}}}(X))_{s}$ if and only if there is an $n\in \mathbb{N}-1$, a word $\mathbf{w}\in S^{n+1}$, and a family $(T^{(1)}_{j})_{j\in n}$ such that $w_{0}=t$, $w_{n}=s$, and, for every $j\in n$, $T^{(1)}_{j}\in \mathrm{Etl}_{w_{j}}(\mathbf{T}_{\Sigma^{\boldsymbol{\mathcal{A}}}}(X))_{w_{j+1}}$ and 
$T^{(1)}=T^{(1)}_{n-1}\circ \cdots \circ T^{(1)}_{0}$. We will sometimes refer to the composition $T^{(1)}_{n-2}\circ \cdots \circ T^{(1)}_{0}$ as the \emph{first-order maximal prefix} of $T^{(1)}$, denoted by $T^{(1)'}$, and we will represent $T^{(1)}$ as $T^{(1)}_{n-1}\circ T^{(1)'}$ or under the form: 
\begin{multline*}
T^{(1)}=\sigma^{\mathbf{T}_{\Sigma^{\boldsymbol{\mathcal{A}}}}(X)}\left(
\mathrm{CH}^{(1)}_{s_{0}}\left(\mathfrak{P}_{0}\right),
\cdots,
\mathrm{CH}^{(1)}_{s_{k-1}}\left(\mathfrak{P}_{k-1}\right),
\right.
\\
\left.
T^{(1)'},
\mathrm{CH}^{(1)}_{s_{k+1}}\left(\mathfrak{P}_{k+1}\right),
\cdots,
\mathrm{CH}^{(1)}_{s_{\bb{\mathbf{s}}-1}}\left(\mathfrak{P}_{\bb{\mathbf{s}}-1}\right)
\right);
\tag{1}
\end{multline*}
\begin{align*}
T^{(1)}&=T^{(1)'}\circ^{0\mathbf{T}_{\Sigma^{\boldsymbol{\mathcal{A}}}}(X)}_{s}\mathrm{CH}^{(1)}_{s}\left(\mathfrak{P}\right);
\\
T^{(1)}&=\mathrm{CH}^{(1)}_{s}\left(\mathfrak{P}\right)\circ^{0\mathbf{T}_{\Sigma^{\boldsymbol{\mathcal{A}}}}(X)}_{s}T^{(1)'}.
\tag{2}
\end{align*}
The underlined space notation can be extended to first-order translations as well. We will say that $T^{(1)}$ is a \emph{first-order translation of type $\sigma$} or a \emph{first-order translation of type $\circ^{0}_{s}$} if $\mathrm{T}^{(1)}_{n-1}$ is a first-order elementary translation of type $\sigma$ or a first-order elementary translation of type $\circ^{0}_{s}$, respectively. We will call $n$ the \emph{height of $T^{(1)}$} and we will denote this fact by $\bb{T^{(1)}}=n$. In this regard, first-order elementary translation have height $1$ and, if $T^{(1)}$ is a first-order translation of height $n$, i.e., $\bb{T^{(1)}}=n$, then its first-order maximal prefix has height $n-1$, i.e., $\bb{T^{(1)'}}=n-1$. We will call the elements of $\mathrm{Tl}_{t}(\mathbf{T}_{\Sigma^{\boldsymbol{\mathcal{A}}}}(X))_{s}$ the \emph{first-order $t$-translations of sort $s$ for $\mathbf{T}_{\Sigma^{\boldsymbol{\mathcal{A}}}}(X)$}. Besides, for every $t\in S$, the mapping $\mathrm{id}^{\mathrm{T}_{\Sigma^{\boldsymbol{\mathcal{A}}}}(X)_{t}}$ will be viewed as an element of $\mathrm{Tl}_{t}(\mathbf{T}_{\Sigma^{\boldsymbol{\mathcal{A}}}}(X))_{t}$. The identity first-order translation has no associated type and we will consider that it has height $0$, i.e., $\bb{\mathrm{id}^{\mathrm{T}_{\Sigma^{\boldsymbol{\mathcal{A}}}}(X)_{t}}}=0$.
Moreover, since we consider that, for every $t\in S$, the identity mapping $\mathrm{id}^{\mathrm{T}_{\Sigma^{\boldsymbol{\mathcal{A}}}}(X)_{t}}$ is a traslation, we agree that a first-order elementary translations $T^{(1)}\colon \mathrm{T}_{\Sigma^{\boldsymbol{\mathcal{A}}}}(X)_{t}\mor
\mathrm{T}_{\Sigma^{\boldsymbol{\mathcal{A}}}}(X)_{s}$, of sort $s$ for $\mathbf{T}_{\Sigma^{\boldsymbol{\mathcal{A}}}}(X)$, has as maximal prefix the identity at the domain.
\end{restatable}

\begin{restatable}{lemma}{LUTransWD}
\label{LUTransWD} 
Let $s, t$ be sorts in $S$, $T^{(1)}$ a first-order translation in 
$\mathrm{Tl}_{t}(\mathrm{T}_{\Sigma^{\boldsymbol{\mathcal{A}}}}(X))_{s}$ and $M$, $N$ path terms in
$\mathrm{PT}_{\boldsymbol{\mathcal{A}}, t}$ such that 
 $$
\left(
\mathrm{ip}^{(1,X)@}_{t}\left(
M
\right), 
\mathrm{ip}^{(1,X)@}_{t}\left(
N
\right)\right)
\in\mathrm{Ker}\left(
\mathrm{sc}^{(0,1)}
\right)_{t}
\cap\mathrm{Ker}\left(
\mathrm{tg}^{(0,1)}
\right)_{t}.
$$
Then the following properties hold 
\begin{itemize}
\item[(i)]
$T^{(1)}(M)$ is a path term in $\mathrm{PT}_{\boldsymbol{\mathcal{A}}, s}$ if, and only if, $T^{(1)}(N)$ is a path term in  $\mathrm{PT}_{\boldsymbol{\mathcal{A}}, s}$;
\item[(ii)] If either $T^{(1)}(M)$ or $T^{(1)}(N)$ is a path term in $\mathrm{PT}_{\boldsymbol{\mathcal{A}}, s}$, then 
\[
\left(
\mathrm{ip}^{(1,X)@}_{s}\left(
T^{(1)}(M)\right),
\mathrm{ip}^{(1,X)@}_{s}\left(
T^{(1)}(N)\right)
\right)
\in\mathrm{Ker}\left(
\mathrm{sc}^{(0,1)}
\right)_{s}
\cap\mathrm{Ker}\left(
\mathrm{tg}^{(0,1)}
\right)_{s}.
\]
\end{itemize}
\end{restatable}
\begin{proof}
We will prove the statement by induction on the height of the first-order translation. 

\textsf{Base case of the induction.}

The statement follows easily when $s=t$ and $T^{(1)}=\mathrm{id}^{\mathrm{T}_{\Sigma^{\boldsymbol{\mathcal{A}}}}(X)_{s}}$.

We consider the different possibilities for the case of $T^{(1)}$ being a first-order elementary translation according to Definition~\ref{DUETrans}.

Assume that there is a word $\mathbf{s}\in S^{\star}-\{\lambda\}$, an index $k\in\bb{\mathbf{s}}$, an operation symbol $\sigma\in \Sigma_{\mathbf{s},s}$, a family of paths $(\mathfrak{P}_{j})_{j\in k}\in \prod_{j\in k}\mathrm{Pth}_{\boldsymbol{\mathcal{A}},s_{j}}$ and a family of paths $(\mathfrak{P}_{l})_{l\in \bb{\mathbf{s}}-(k+1)}\in \prod_{l\in \bb{\mathbf{s}}-(k+1)}\mathrm{Pth}_{\boldsymbol{\mathcal{A}},s_{l}}$ such that $s_{k}=t$ and
\begin{multline*}
T^{(1)}=\sigma^{\mathbf{T}_{\Sigma^{\boldsymbol{\mathcal{A}}}}(X)}\left(
\mathrm{CH}^{(1)}_{s_{0}}\left(\mathfrak{P}_{0}\right),
\cdots,
\mathrm{CH}^{(1)}_{s_{k-1}}\left(\mathfrak{P}_{k-1}\right),
\right.
\\
\left.
\underline{\quad},
\mathrm{CH}^{(1)}_{s_{k+1}}\left(\mathfrak{P}_{k+1}\right),
\cdots,
\mathrm{CH}^{(1)}_{s_{\bb{\mathbf{s}}-1}}\left(\mathfrak{P}_{\bb{\mathbf{s}}-1}\right)
\right).
\end{multline*}

Note that, for this case, the terms $T^{(1)}(M)$ and $T^{(1)}(N)$ are always path terms. Therefore, the first condition of the proposition holds.

Now, regarding the $(0,1)$-source, the following chain of equalities holds
\begin{flushleft}
$\mathrm{sc}^{(0,1)}_{s}\left(\mathrm{ip}^{(1,X)@}_{s}\left(
T^{(1)}(M)\right)\right)$
\allowdisplaybreaks
\begin{align*}
\qquad&=
\mathrm{sc}^{(0,1)}_{s}\left(
\mathrm{ip}^{(1,X)@}_{s}\left(
\sigma^{\mathbf{T}_{\Sigma^{\boldsymbol{\mathcal{A}}}}(X)}\left(
\mathrm{CH}^{(1)}_{s_{0}}\left(\mathfrak{P}_{0}\right),
\cdots,
\mathrm{CH}^{(1)}_{s_{k-1}}\left(\mathfrak{P}_{k-1}\right),
\right.\right.\right.
\\&\qquad\qquad\qquad\qquad\qquad
\left.\left.\left.
M,
\mathrm{CH}^{(1)}_{s_{k+1}}\left(\mathfrak{P}_{k+1}\right),
\cdots,
\mathrm{CH}^{(1)}_{s_{\bb{\mathbf{s}}-1}}\left(\mathfrak{P}_{\bb{\mathbf{s}}-1}\right)
\right)
\right)\right)
\tag{1}
\\&=
\mathrm{sc}^{(0,1)}_{s}\left(
\sigma^{\mathbf{Pth}_{\boldsymbol{\mathcal{A}}}}\left(
\mathrm{ip}^{(1,X)@}_{s_{0}}\left(
\mathrm{CH}^{(1)}_{s_{0}}\left(\mathfrak{P}_{0}\right)\right),
\cdots,
\mathrm{ip}^{(1,X)@}_{s_{k-1}}\left(
\mathrm{CH}^{(1)}_{s_{k-1}}\left(\mathfrak{P}_{k-1}\right)\right),
\right.\right.
\\&\qquad\qquad\qquad\qquad
\mathrm{ip}^{(1,X)@}_{t}\left(
M\right),
\mathrm{ip}^{(1,X)@}_{s_{k+1}}\left(
\mathrm{CH}^{(1)}_{s_{k+1}}\left(\mathfrak{P}_{k+1}\right)\right),
\\&\qquad\qquad\qquad\qquad\qquad\qquad\qquad\quad
\left.\left.
\cdots,
\mathrm{ip}^{(1,X)@}_{s_{\bb{\mathbf{s}}-1}}\left(
\mathrm{CH}^{(1)}_{s_{\bb{\mathbf{s}}-1}}\left(\mathfrak{P}_{\bb{\mathbf{s}}-1}\right)\right)
\right)
\right)
\tag{2}
\\&=
\sigma^{\mathbf{T}_{\Sigma}(X)}\left(
\mathrm{sc}^{(0,1)}_{s_{0}}\left(
\mathrm{ip}^{(1,X)@}_{s_{0}}\left(
\mathrm{CH}^{(1)}_{s_{0}}\left(\mathfrak{P}_{0}\right)\right)\right),
\cdots,
\right.
\\&\qquad\qquad
\mathrm{sc}^{(0,1)}_{s_{k-1}}\left(
\mathrm{ip}^{(1,X)@}_{s_{k-1}}\left(
\mathrm{CH}^{(1)}_{s_{k-1}}\left(\mathfrak{P}_{k-1}\right)\right)\right),
\\&\qquad\qquad\qquad
\mathrm{sc}^{(0,1)}_{t}\left(
\mathrm{ip}^{(1,X)@}_{t}\left(
M\right)\right),
\mathrm{sc}^{(0,1)}_{s_{k+1}}\left(
\mathrm{ip}^{(1,X)@}_{s_{k+1}}\left(
\mathrm{CH}^{(1)}_{s_{k+1}}\left(\mathfrak{P}_{k+1}\right)\right)\right),
\\&\qquad\qquad\qquad\qquad\qquad\qquad
\left.
\cdots,
\mathrm{sc}^{(0,1)}_{s_{\bb{\mathbf{s}}-1}}\left(
\mathrm{ip}^{(1,X)@}_{s_{\bb{\mathbf{s}}-1}}\left(
\mathrm{CH}^{(1)}_{s_{\bb{\mathbf{s}}-1}}\left(\mathfrak{P}_{\bb{\mathbf{s}}-1}\right)\right)
\right)
\right)
\tag{3}
\\&=
\sigma^{\mathbf{T}_{\Sigma}(X)}\left(
\mathrm{sc}^{(0,1)}_{s_{0}}\left(
\mathrm{ip}^{(1,X)@}_{s_{0}}\left(
\mathrm{CH}^{(1)}_{s_{0}}\left(\mathfrak{P}_{0}\right)\right)\right),
\cdots,
\right.
\\&\qquad\qquad
\mathrm{sc}^{(0,1)}_{s_{k-1}}\left(
\mathrm{ip}^{(1,X)@}_{s_{k-1}}\left(
\mathrm{CH}^{(1)}_{s_{k-1}}\left(\mathfrak{P}_{k-1}\right)\right)\right),
\\&\qquad\qquad\qquad
\mathrm{sc}^{(0,1)}_{t}\left(
\mathrm{ip}^{(1,X)@}_{t}\left(
N\right)\right),
\mathrm{sc}^{(0,1)}_{s_{k+1}}\left(
\mathrm{ip}^{(1,X)@}_{s_{k+1}}\left(
\mathrm{CH}^{(1)}_{s_{k+1}}\left(\mathfrak{P}_{k+1}\right)\right)\right),
\\&\qquad\qquad\qquad\qquad\qquad\qquad
\left.
\cdots,
\mathrm{sc}^{(0,1)}_{s_{\bb{\mathbf{s}}-1}}\left(
\mathrm{ip}^{(1,X)@}_{s_{\bb{\mathbf{s}}-1}}\left(
\mathrm{CH}^{(1)}_{s_{\bb{\mathbf{s}}-1}}\left(\mathfrak{P}_{\bb{\mathbf{s}}-1}\right)\right)
\right)
\right)
\tag{4}
\\&=
\mathrm{sc}^{(0,1)}_{s}\left(
\sigma^{\mathbf{Pth}_{\boldsymbol{\mathcal{A}}}}\left(
\mathrm{ip}^{(1,X)@}_{s_{0}}\left(
\mathrm{CH}^{(1)}_{s_{0}}\left(\mathfrak{P}_{0}\right)\right),
\cdots,
\mathrm{ip}^{(1,X)@}_{s_{k-1}}\left(
\mathrm{CH}^{(1)}_{s_{k-1}}\left(\mathfrak{P}_{k-1}\right)\right),
\right.\right.
\\&\qquad\qquad\qquad\qquad
\mathrm{ip}^{(1,X)@}_{t}\left(
N\right),
\mathrm{ip}^{(1,X)@}_{s_{k+1}}\left(
\mathrm{CH}^{(1)}_{s_{k+1}}\left(\mathfrak{P}_{k+1}\right)\right),
\\&\qquad\qquad\qquad\qquad\qquad\qquad\qquad\quad
\left.\left.
\cdots,
\mathrm{ip}^{(1,X)@}_{s_{\bb{\mathbf{s}}-1}}\left(
\mathrm{CH}^{(1)}_{s_{\bb{\mathbf{s}}-1}}\left(\mathfrak{P}_{\bb{\mathbf{s}}-1}\right)\right)
\right)
\right)
\tag{5}
\\&
\mathrm{sc}^{(0,1)}_{s}\left(
\mathrm{ip}^{(1,X)@}_{s}\left(
\sigma^{\mathbf{T}_{\Sigma^{\boldsymbol{\mathcal{A}}}}(X)}\left(
\mathrm{CH}^{(1)}_{s_{0}}\left(\mathfrak{P}_{0}\right),
\cdots,
\mathrm{CH}^{(1)}_{s_{k-1}}\left(\mathfrak{P}_{k-1}\right),
\right.\right.\right.
\\&\qquad\qquad\qquad\qquad\qquad
\left.\left.\left.
N,
\mathrm{CH}^{(1)}_{s_{k+1}}\left(\mathfrak{P}_{k+1}\right),
\cdots,
\mathrm{CH}^{(1)}_{s_{\bb{\mathbf{s}}-1}}\left(\mathfrak{P}_{\bb{\mathbf{s}}-1}\right)
\right)
\right)\right)
\tag{6}
\\&=\mathrm{sc}^{(0,1)}_{s}\left(\mathrm{ip}^{(1,X)@}_{s}\left(
T^{(1)}(N)\right)\right).
\tag{7}
\end{align*}
\end{flushleft}

In the just stated chain of equalities, the first equality unravels the description of the first-order elementary translation $T^{(1)}$; the second equality follows from the fact that $\mathrm{ip}^{(1,X)@}$ is a $\Sigma^{\boldsymbol{\mathcal{A}}}$-homomorphism, according to Definition~\ref{DIp}; the third equality follows from the fact that $\mathrm{sc}^{(0,1)}$ is a $\Sigma$-homomorphism, according to Proposition~\ref{PHom}; the fourth equality follows from the fact that, by hypothesis, 
$
\mathrm{sc}^{(0,1)}_{t}(
\mathrm{ip}^{(1,X)@}_{t}(
M))
=
\mathrm{sc}^{(0,1)}_{t}(
\mathrm{ip}^{(1,X)@}_{t}(
N))$;
the fifth equality follows from the fact that $\mathrm{sc}^{(0,1)}$ is a $\Sigma$-homomorphism, according to Proposition~\ref{PHom}; the sixth equality follows from the fact that $\mathrm{ip}^{(1,X)@}$ is a $\Sigma^{\boldsymbol{\mathcal{A}}}$-homomorphism, according to Definition~\ref{DIp}; finally, the last equality recovers  the description of the first-order elementary translation $T^{(1)}$.

The case of the $(0,1)$-target follows by a similar argument.

Now, regarding the other options for the first-order elementary translation $T^{(1)}$, it could be the case that $s=t$ and there is a path $\mathfrak{P}\in\mathrm{Pth}_{\boldsymbol{\mathcal{A}},s}$ such that
\[
T^{(1)}=\underline{\quad}\circ^{0\mathbf{T}_{\Sigma^{\boldsymbol{\mathcal{A}}}}(X)}_{s}\mathrm{CH}^{(1)}_{s}\left(\mathfrak{P}\right).
\]

For this case, the following chain of equivalences hold
\begin{flushleft}
$T^{(1)}(M)$ is a path term
\allowdisplaybreaks
\begin{align*}
\Leftrightarrow\quad & M\circ^{0\mathbf{T}_{\Sigma^{\boldsymbol{\mathcal{A}}}}(X)}_{s}\mathrm{CH}^{(1)}_{s}\left(\mathfrak{P}\right)\mbox{ is a path term}
\tag{1}
\\
\Leftrightarrow\quad &
\mathrm{ip}^{(1,X)@}_{s}\left(
M\circ^{0\mathbf{T}_{\Sigma^{\boldsymbol{\mathcal{A}}}}(X)}_{s}\mathrm{CH}^{(1)}_{s}\left(\mathfrak{P}\right)
\right)
\mbox{ is a path}
\tag{2}
\\
\Leftrightarrow\quad &
\mathrm{ip}^{(1,X)@}_{s}\left(
M\right)\circ^{0\mathbf{Pth}_{\boldsymbol{\mathcal{A}}}}_{s}
\mathrm{ip}^{(1,X)@}_{s}\left(
\mathrm{CH}^{(1)}_{s}\left(\mathfrak{P}\right)
\right)
\mbox{ is a path}
\tag{3}
\\
\Leftrightarrow\quad &
\mathrm{sc}^{(0,1)}_{s}\left(
\mathrm{ip}^{(1,X)@}_{s}\left(
M\right)\right)
=
\mathrm{tg}^{(0,1)}_{s}\left(
\mathrm{ip}^{(1,X)@}_{s}\left(
\mathrm{CH}^{(1)}_{s}\left(\mathfrak{P}\right)
\right)\right)
\tag{4}
\\
\Leftrightarrow\quad &
\mathrm{sc}^{(0,1)}_{s}\left(
\mathrm{ip}^{(1,X)@}_{s}\left(
N\right)\right)
=
\mathrm{tg}^{(0,1)}_{s}\left(
\mathrm{ip}^{(1,X)@}_{s}\left(
\mathrm{CH}^{(1)}_{s}\left(\mathfrak{P}\right)
\right)\right)
\tag{5}
\\
\Leftrightarrow\quad &
\mathrm{ip}^{(1,X)@}_{s}\left(
N\right)\circ^{0\mathbf{Pth}_{\boldsymbol{\mathcal{A}}}}_{s}
\mathrm{ip}^{(1,X)@}_{s}\left(
\mathrm{CH}^{(1)}_{s}\left(\mathfrak{P}\right)
\right)
\mbox{ is a path}
\tag{6}
\\
\Leftrightarrow\quad &
\mathrm{ip}^{(1,X)@}_{s}\left(
N\circ^{0\mathbf{T}_{\Sigma^{\boldsymbol{\mathcal{A}}}}(X)}_{s}\mathrm{CH}^{(1)}_{s}\left(\mathfrak{P}\right)
\right)
\mbox{ is a path}
\tag{7}
\\
\Leftrightarrow\quad & N\circ^{0\mathbf{T}_{\Sigma^{\boldsymbol{\mathcal{A}}}}(X)}_{s}\mathrm{CH}^{(1)}_{s}\left(\mathfrak{P}\right)\mbox{ is a path term}
\tag{8}
\\
\Leftrightarrow\quad &
T^{(1)}(N) \mbox{ is a path term}.
\tag{9}
\end{align*}
\end{flushleft}

In the just stated chain of equivalences, the first equivalence follows from the description of the first-order elementary translation $T^{(1)}$; the second equivalence follows from Proposition~\ref{PPT}; 	the third equivalence follows from the fact that $\mathrm{ip}^{(1,X)@}$ is a $\Sigma^{\boldsymbol{\mathcal{A}}}$-homomorphism, according to Definition~\ref{DIp}; the fourth equivalence follows from the definition of the $0$-composition operation in $\mathbf{Pth}_{\boldsymbol{\mathcal{A}}}$ according to Proposition~\ref{PPthCatAlg}; the  fifth equivalence follows from the fact that, by hypothesis, $
\mathrm{sc}^{(0,1)}_{t}(
\mathrm{ip}^{(1,X)@}_{t}(
M))
=
\mathrm{sc}^{(0,1)}_{t}(
\mathrm{ip}^{(1,X)@}_{t}(
N))$; the sixth equivalence follows from the definition of the $0$-composition operation in $\mathbf{Pth}_{\boldsymbol{\mathcal{A}}}$ according to Proposition~\ref{PPthCatAlg}; the seventh equivalence follows from the fact that $\mathrm{ip}^{(1,X)@}$ is a $\Sigma^{\boldsymbol{\mathcal{A}}}$-homomorphism, according to Definition~\ref{DIp}; the eighth equivalence follows from Proposition~\ref{PPT}; finally, the last equivalence follows from the description of the first-order elementary translation $T^{(1)}$.

Assume that either $T^{(1)}(M)$ or $T^{(1)}(N)$ is a path term.
Regarding the $(0,1)$-source, the following chain of equalities holds
\begin{flushleft}
$\mathrm{sc}^{(0,1)}_{s}\left(\mathrm{ip}^{(1,X)@}_{s}\left(
T^{(1)}(M)\right)\right)$
\allowdisplaybreaks
\begin{align*}
\qquad&=\mathrm{sc}^{(0,1)}_{s}\left(\mathrm{ip}^{(1,X)@}_{s}\left(
M\circ^{0\mathbf{T}_{\Sigma^{\boldsymbol{\mathcal{A}}}}(X)}_{s}\mathrm{CH}^{(1)}_{s}\left(\mathfrak{P}\right)\right)\right)
\tag{1}
\\&=
\mathrm{sc}^{(0,1)}_{s}\left(
\mathrm{ip}^{(1,X)@}_{s}\left(
M\right)
\circ^{0\mathbf{Pth}_{\boldsymbol{\mathcal{A}}}}_{s}
\mathrm{ip}^{(1,X)@}_{s}\left(
\mathrm{CH}^{(1)}_{s}\left(
\mathfrak{P}\right)
\right)\right)
\tag{2}
\\&=
\mathrm{sc}^{(0,1)}_{s}\left(
\mathrm{ip}^{(1,X)@}_{s}\left(
\mathrm{CH}^{(1)}_{s}\left(
\mathfrak{P}\right)
\right)\right)
\tag{3}
\\&=
\mathrm{sc}^{(0,1)}_{s}\left(
\mathrm{ip}^{(1,X)@}_{s}\left(
N\right)
\circ^{0\mathbf{Pth}_{\boldsymbol{\mathcal{A}}}}_{s}
\mathrm{ip}^{(1,X)@}_{s}\left(
\mathrm{CH}^{(1)}_{s}\left(
\mathfrak{P}\right)
\right)\right)
\tag{4}
\\&=\mathrm{sc}^{(0,1)}_{s}\left(\mathrm{ip}^{(1,X)@}_{s}\left(
N\circ^{0\mathbf{T}_{\Sigma^{\boldsymbol{\mathcal{A}}}}(X)}_{s}\mathrm{CH}^{(1)}_{s}\left(\mathfrak{P}\right)\right)\right)
\tag{5}
\\&=
\mathrm{sc}^{(0,1)}_{s}\left(\mathrm{ip}^{(1,X)@}_{s}\left(
T^{(1)}(N)\right)\right).
\tag{6}
\end{align*}
\end{flushleft}

In the just stated chain of equalities, the first equality follows from the description of the first-order elementary translation $T^{(1)}$; the second equality follows from the fact that $\mathrm{ip}^{(1,X)@}$ is a $\Sigma^{\boldsymbol{\mathcal{A}}}$-homomorphism, according to Definition~\ref{DIp}; the third equality follows from the definition of the $0$-composition operation in $\mathbf{Pth}_{\boldsymbol{\mathcal{A}}}$ according to Proposition~\ref{PPthCatAlg}; the fourth equality follows from the definition of the $0$-composition operation in $\mathbf{Pth}_{\boldsymbol{\mathcal{A}}}$ according to Proposition~\ref{PPthCatAlg}; the fifth equality follows from the fact that $\mathrm{ip}^{(1,X)@}$ is a $\Sigma^{\boldsymbol{\mathcal{A}}}$-homomorphism, according to Definition~\ref{DIp}; finally, the last equivalence follows from the description of the first-order elementary translation $T^{(1)}$.

Now, regarding the $(0,1)$-target, the following chain of equalities holds
\begin{flushleft}
$\mathrm{tg}^{(0,1)}_{s}\left(\mathrm{ip}^{(1,X)@}_{s}\left(
T^{(1)}(M)\right)\right)$
\allowdisplaybreaks
\begin{align*}
\qquad&=\mathrm{tg}^{(0,1)}_{s}\left(\mathrm{ip}^{(1,X)@}_{s}\left(
M\circ^{0\mathbf{T}_{\Sigma^{\boldsymbol{\mathcal{A}}}}(X)}_{s}\mathrm{CH}^{(1)}_{s}\left(\mathfrak{P}\right)\right)\right)
\tag{1}
\\&=
\mathrm{tg}^{(0,1)}_{s}\left(
\mathrm{ip}^{(1,X)@}_{s}\left(
M\right)
\circ^{0\mathbf{Pth}_{\boldsymbol{\mathcal{A}}}}_{s}
\mathrm{ip}^{(1,X)@}_{s}\left(
\mathrm{CH}^{(1)}_{s}\left(
\mathfrak{P}\right)
\right)\right)
\tag{2}
\\&=
\mathrm{tg}^{(0,1)}_{s}\left(
\mathrm{ip}^{(1,X)@}_{s}\left(
M\right)
\right)
\tag{3}
\\&=
\mathrm{tg}^{(0,1)}_{s}\left(
\mathrm{ip}^{(1,X)@}_{s}\left(
N\right)
\right)
\tag{4}
\\&=
\mathrm{tg}^{(0,1)}_{s}\left(
\mathrm{ip}^{(1,X)@}_{s}\left(
N\right)
\circ^{0\mathbf{Pth}_{\boldsymbol{\mathcal{A}}}}_{s}
\mathrm{ip}^{(1,X)@}_{s}\left(
\mathrm{CH}^{(1)}_{s}\left(
\mathfrak{P}\right)
\right)\right)
\tag{5}
\\&=\mathrm{tg}^{(0,1)}_{s}\left(\mathrm{ip}^{(1,X)@}_{s}\left(
N\circ^{0\mathbf{T}_{\Sigma^{\boldsymbol{\mathcal{A}}}}(X)}_{s}\mathrm{CH}^{(1)}_{s}\left(\mathfrak{P}\right)\right)\right)
\tag{6}
\\&=
\mathrm{tg}^{(0,1)}_{s}\left(\mathrm{ip}^{(1,X)@}_{s}\left(
T^{(1)}(N)\right)\right).
\tag{7}
\end{align*}
\end{flushleft}

In the just stated chain of equalities, the first equality follows from the description of the first-order elementary translation $T^{(1)}$; the second equality follows from the fact that $\mathrm{ip}^{(1,X)@}$ is a $\Sigma^{\boldsymbol{\mathcal{A}}}$-homomorphism, according to Definition~\ref{DIp}; the third equality follows from the definition of the $0$-composition operation in $\mathbf{Pth}_{\boldsymbol{\mathcal{A}}}$ according to Proposition~\ref{PPthCatAlg}; the fourth
equality follows from the fact that, by hypothesis, $
\mathrm{tg}^{(0,1)}_{t}(
\mathrm{ip}^{(1,X)@}_{t}(
M))
=
\mathrm{tg}^{(0,1)}_{t}(
\mathrm{ip}^{(1,X)@}_{t}(
N))$; the fifth equality follows from the definition of the $0$-composition operation in $\mathbf{Pth}_{\boldsymbol{\mathcal{A}}}$ according to Proposition~\ref{PPthCatAlg}; the sixth equality follows from the fact that $\mathrm{ip}^{(1,X)@}$ is a $\Sigma^{\boldsymbol{\mathcal{A}}}$-homomorphism, according to Definition~\ref{DIp}; finally, the last equivalence follows from the description of the first-order elementary translation $T^{(1)}$.

The remaining case, i.e., the case for which $s=t$ and there is a path $\mathfrak{P}\in\mathrm{Pth}_{\boldsymbol{\mathcal{A}},s}$ such that
\[
T^{(1)}=\mathrm{CH}^{(1)}_{s}\left(\mathfrak{P}\right)\circ^{0\mathbf{T}_{\Sigma^{\boldsymbol{\mathcal{A}}}}(X)}_{s}\underline{\quad},
\]
follows by a similar argument to the proof presented above.

This complete the base case.

\textsf{Inductive step of the induction.}

Let us suppose that the statement holds for first-order translations of height $m\in \mathbb{N}-\{0\}$. That is, if $T^{(1)}$ is a first-order translation in 
$\mathrm{Tl}_{t}(\mathrm{T}_{\Sigma^{\boldsymbol{\mathcal{A}}}}(X))_{s}$ of height $m$ and $M$ and $N$ are path terms in
$\mathrm{PT}_{\boldsymbol{\mathcal{A}}, t}$ satisfying that 
 $$
\left(
\mathrm{ip}^{(1,X)@}_{t}\left(
M
\right), 
\mathrm{ip}^{(1,X)@}_{t}\left(
N
\right)\right)
\in\mathrm{Ker}\left(
\mathrm{sc}^{(0,1)}
\right)_{t}
\cap\mathrm{Ker}\left(
\mathrm{tg}^{(0,1)}
\right)_{t},
$$
then the following properties hold 
\begin{itemize}
\item[(i)]
$T^{(1)}(M)$ is a path term in $\mathrm{PT}_{\boldsymbol{\mathcal{A}}, s}$ if, and only if, $T^{(1)}(N)$ is a path term in  $\mathrm{PT}_{\boldsymbol{\mathcal{A}}, s}$;
\item[(ii)] If either $T^{(1)}(M)$ or $T^{(1)}(N)$ is a path term in $\mathrm{PT}_{\boldsymbol{\mathcal{A}}, s}$, then 
\[
\left(
\mathrm{ip}^{(1,X)@}_{s}\left(
T^{(1)}(M)\right),
\mathrm{ip}^{(1,X)@}_{s}\left(
T^{(1)}(N)\right)
\right)
\in\mathrm{Ker}\left(
\mathrm{sc}^{(0,1)}
\right)_{s}
\cap\mathrm{Ker}\left(
\mathrm{tg}^{(0,1)}
\right)_{s}.
\]
\end{itemize}

Now let $T^{(1)}$ be a first-order translation in 
$\mathrm{Tl}_{t}(\mathrm{T}_{\Sigma^{\boldsymbol{\mathcal{A}}}}(X))_{s}$ of height $m+1$ and let $M$ and $N$ be path terms in
$\mathrm{PT}_{\boldsymbol{\mathcal{A}}, t}$ satisfying that 
 $$
\left(
\mathrm{ip}^{(1,X)@}_{t}\left(
M
\right), 
\mathrm{ip}^{(1,X)@}_{t}\left(
N
\right)\right)
\in\mathrm{Ker}\left(
\mathrm{sc}^{(0,1)}
\right)_{t}
\cap\mathrm{Ker}\left(
\mathrm{tg}^{(0,1)}
\right)_{t}.
$$

We consider the different possibilities for the case of $T^{(1)}$ being a first-order  translation according to Definition~\ref{DUTrans}.

Assume that there is a word $\mathbf{s}\in S^{\star}-\{\lambda\}$, an index $k\in\bb{\mathbf{s}}$, an operation symbol $\sigma\in \Sigma_{\mathbf{s},s}$, a family of paths $(\mathfrak{P}_{j})_{j\in k}\in \prod_{j\in k}\mathrm{Pth}_{\boldsymbol{\mathcal{A}},s_{j}}$, a family of paths $(\mathfrak{P}_{l})_{l\in \bb{\mathbf{s}}-(k+1)}\in \prod_{l\in \bb{\mathbf{s}}-(k+1)}\mathrm{Pth}_{\boldsymbol{\mathcal{A}},s_{l}}$ 
and a first-order translation $T^{(1)'}$ in $\mathrm{Tl}_{t}(\mathrm{T}_{\Sigma^{\boldsymbol{\mathcal{A}}}}(X))_{s_{k}}$ of height $m$  such that 
\begin{multline*}
T^{(1)}=\sigma^{\mathbf{T}_{\Sigma^{\boldsymbol{\mathcal{A}}}}(X)}\left(
\mathrm{CH}^{(1)}_{s_{0}}\left(\mathfrak{P}_{0}\right),
\cdots,
\mathrm{CH}^{(1)}_{s_{k-1}}\left(\mathfrak{P}_{k-1}\right),
\right.
\\
\left.
T^{(1)'},
\mathrm{CH}^{(1)}_{s_{k+1}}\left(\mathfrak{P}_{k+1}\right),
\cdots,
\mathrm{CH}^{(1)}_{s_{\bb{\mathbf{s}}-1}}\left(\mathfrak{P}_{\bb{\mathbf{s}}-1}\right)
\right).
\end{multline*}

Note that, for this case, the following equivalences hold
\begin{flushleft}
$T^{(1)}(M)$ is a path term
\allowdisplaybreaks
\begin{align*}
\Leftrightarrow\qquad&
\sigma^{\mathbf{T}_{\Sigma^{\boldsymbol{\mathcal{A}}}}(X)}\left(
\mathrm{CH}^{(1)}_{s_{0}}\left(\mathfrak{P}_{0}\right),
\cdots,
\mathrm{CH}^{(1)}_{s_{k-1}}\left(\mathfrak{P}_{k-1}\right),
T^{(1)'}(M),
\right.
\\&\qquad\qquad\qquad
\left.
\mathrm{CH}^{(1)}_{s_{k+1}}\left(\mathfrak{P}_{k+1}\right),
\cdots,
\mathrm{CH}^{(1)}_{s_{\bb{\mathbf{s}}-1}}\left(\mathfrak{P}_{\bb{\mathbf{s}}-1}\right)
\right)
\mbox{ is a path term}
\tag{1}
\\
\Leftrightarrow\qquad&
\mathrm{ip}^{(1,X)@}_{s}\left(
\sigma^{\mathbf{T}_{\Sigma^{\boldsymbol{\mathcal{A}}}}(X)}\left(
\mathrm{CH}^{(1)}_{s_{0}}\left(\mathfrak{P}_{0}\right),
\cdots,
\mathrm{CH}^{(1)}_{s_{k-1}}\left(\mathfrak{P}_{k-1}\right),
T^{(1)'}(M),
\right.\right.
\\&\qquad\qquad\qquad\qquad
\left.\left.
\mathrm{CH}^{(1)}_{s_{k+1}}\left(\mathfrak{P}_{k+1}\right),
\cdots,
\mathrm{CH}^{(1)}_{s_{\bb{\mathbf{s}}-1}}\left(\mathfrak{P}_{\bb{\mathbf{s}}-1}\right)
\right)\right)
\mbox{ is a path}
\tag{2}
\\
\Leftrightarrow\qquad&
\sigma^{\mathbf{Pth}_{\boldsymbol{\mathcal{A}}}}\left(
\mathrm{ip}^{(1,X)@}_{s_{0}}\left(
\mathrm{CH}^{(1)}_{s_{0}}\left(\mathfrak{P}_{0}\right)\right),
\cdots,
\mathrm{ip}^{(1,X)@}_{s_{k-1}}\left(
\mathrm{CH}^{(1)}_{s_{k-1}}\left(\mathfrak{P}_{k-1}\right)\right),
\right.
\\&\qquad\qquad\qquad
\mathrm{ip}^{(1,X)@}_{s_{k}}\left(
T^{(1)'}(M)\right),
\mathrm{ip}^{(1,X)@}_{s_{k+1}}\left(
\mathrm{CH}^{(1)}_{s_{k+1}}\left(\mathfrak{P}_{k+1}\right)\right),
\cdots,
\\&\qquad\qquad\qquad\qquad\qquad\qquad
\left.
\mathrm{ip}^{(1,X)@}_{s_{\bb{\mathbf{s}}-1}}\left(
\mathrm{CH}^{(1)}_{s_{\bb{\mathbf{s}}-1}}\left(\mathfrak{P}_{\bb{\mathbf{s}}-1}\right)
\right)\right)
\mbox{ is a path}
\tag{3}
\\
\Leftrightarrow\qquad&
\mathrm{ip}^{(1,X)@}_{s_{k}}\left(
T^{(1)'}(M)\right)
\mbox{ is a path}
\tag{4}
\\
\Leftrightarrow\qquad&
T^{(1)'}(M)
\mbox{ is a path term}
\tag{5}
\\
\Leftrightarrow\qquad&
T^{(1)'}(N)
\mbox{ is a path term}
\tag{6}
\\
\Leftrightarrow\qquad&
\mathrm{ip}^{(1,X)@}_{s_{k}}\left(
T^{(1)'}(N)\right)
\mbox{ is a path}
\tag{7}
\\
\Leftrightarrow\qquad&
\sigma^{\mathbf{Pth}_{\boldsymbol{\mathcal{A}}}}\left(
\mathrm{ip}^{(1,X)@}_{s_{0}}\left(
\mathrm{CH}^{(1)}_{s_{0}}\left(\mathfrak{P}_{0}\right)\right),
\cdots,
\mathrm{ip}^{(1,X)@}_{s_{k-1}}\left(
\mathrm{CH}^{(1)}_{s_{k-1}}\left(\mathfrak{P}_{k-1}\right)\right),
\right.
\\&\qquad\qquad\qquad
\mathrm{ip}^{(1,X)@}_{s_{k}}\left(
T^{(1)'}(N)\right),
\mathrm{ip}^{(1,X)@}_{s_{k+1}}\left(
\mathrm{CH}^{(1)}_{s_{k+1}}\left(\mathfrak{P}_{k+1}\right)\right),
\cdots,
\\&\qquad\qquad\qquad\qquad\qquad\qquad
\left.
\mathrm{ip}^{(1,X)@}_{s_{\bb{\mathbf{s}}-1}}\left(
\mathrm{CH}^{(1)}_{s_{\bb{\mathbf{s}}-1}}\left(\mathfrak{P}_{\bb{\mathbf{s}}-1}\right)
\right)\right)
\mbox{ is a path}
\tag{8}
\\
\Leftrightarrow\qquad&
\mathrm{ip}^{(1,X)@}_{s}\left(
\sigma^{\mathbf{T}_{\Sigma^{\boldsymbol{\mathcal{A}}}}(X)}\left(
\mathrm{CH}^{(1)}_{s_{0}}\left(\mathfrak{P}_{0}\right),
\cdots,
\mathrm{CH}^{(1)}_{s_{k-1}}\left(\mathfrak{P}_{k-1}\right),
T^{(1)'}(N),
\right.\right.
\\&\qquad\qquad\qquad\qquad
\left.\left.
\mathrm{CH}^{(1)}_{s_{k+1}}\left(\mathfrak{P}_{k+1}\right),
\cdots,
\mathrm{CH}^{(1)}_{s_{\bb{\mathbf{s}}-1}}\left(\mathfrak{P}_{\bb{\mathbf{s}}-1}\right)
\right)\right)
\mbox{ is a path}
\tag{9}
\\
\Leftrightarrow\qquad&
\sigma^{\mathbf{T}_{\Sigma^{\boldsymbol{\mathcal{A}}}}(X)}\left(
\mathrm{CH}^{(1)}_{s_{0}}\left(\mathfrak{P}_{0}\right),
\cdots,
\mathrm{CH}^{(1)}_{s_{k-1}}\left(\mathfrak{P}_{k-1}\right),
T^{(1)'}(N),
\right.
\\&\qquad\qquad\qquad
\left.
\mathrm{CH}^{(1)}_{s_{k+1}}\left(\mathfrak{P}_{k+1}\right),
\cdots,
\mathrm{CH}^{(1)}_{s_{\bb{\mathbf{s}}-1}}\left(\mathfrak{P}_{\bb{\mathbf{s}}-1}\right)
\right)
\mbox{ is a path term}
\tag{10}
\\
\Leftrightarrow\qquad&
T^{(1)}(N) \mbox{ is a path term}.
\tag{11}
\end{align*}
\end{flushleft}

In the just stated chain of equivalences, the first equivalence follows from the description of the first-order translation $T^{(1)}$; the second equivalence follows from Proposition~\ref{PPT}; the third equivalence follows from the fact that $\mathrm{ip}^{(1,X)@}$ is a $\Sigma^{\boldsymbol{\mathcal{A}}}$-homomorphism, according to Definition~\ref{DIp}; the fourth equivalence follows from the description of the operation $\sigma$ in the many-sorted partial $\Sigma^{\boldsymbol{\mathcal{A}}}$-algebra $\mathbf{Pth}_{\boldsymbol{\mathcal{A}}}$, according to Proposition~\ref{PPthAlg}; the fifth equivalence follows from Proposition~\ref{PPT}; the sixth equivalence follows by induction; the seventh equivalence follows from Proposition~\ref{PPT}; the eighth equivalence follows from the description of the operation $\sigma$ in the many-sorted partial $\Sigma^{\boldsymbol{\mathcal{A}}}$-algebra $\mathbf{Pth}_{\boldsymbol{\mathcal{A}}}$, according to Proposition~\ref{PPthAlg}; the ninth equivalence follows from the fact that $\mathrm{ip}^{(1,X)@}$ is a $\Sigma^{\boldsymbol{\mathcal{A}}}$-homomorphism, according to Definition~\ref{DIp}; the tenth equivalence follows from Proposition~\ref{PPT}; finally, the last equivalence follows from the description of the first-order translation $T^{(1)}$.

Therefore, the first condition of the proposition holds.

Now, regarding the $(0,1)$-source, the following chain of equalities holds
\begin{flushleft}
$\mathrm{sc}^{(0,1)}_{s}\left(\mathrm{ip}^{(1,X)@}_{s}\left(
T^{(1)}(M)\right)\right)$
\allowdisplaybreaks
\begin{align*}
\qquad&=
\mathrm{sc}^{(0,1)}_{s}\left(
\mathrm{ip}^{(1,X)@}_{s}\left(
\sigma^{\mathbf{T}_{\Sigma^{\boldsymbol{\mathcal{A}}}}(X)}\left(
\mathrm{CH}^{(1)}_{s_{0}}\left(\mathfrak{P}_{0}\right),
\cdots,
\mathrm{CH}^{(1)}_{s_{k-1}}\left(\mathfrak{P}_{k-1}\right),
\right.\right.\right.
\\&\qquad\qquad\qquad\qquad
\left.\left.\left.
T^{(1)'}(M),
\mathrm{CH}^{(1)}_{s_{k+1}}\left(\mathfrak{P}_{k+1}\right),
\cdots,
\mathrm{CH}^{(1)}_{s_{\bb{\mathbf{s}}-1}}\left(\mathfrak{P}_{\bb{\mathbf{s}}-1}\right)
\right)
\right)\right)
\tag{1}
\\&=
\mathrm{sc}^{(0,1)}_{s}\left(
\sigma^{\mathbf{Pth}_{\boldsymbol{\mathcal{A}}}}\left(
\mathrm{ip}^{(1,X)@}_{s_{0}}\left(
\mathrm{CH}^{(1)}_{s_{0}}\left(\mathfrak{P}_{0}\right)\right),
\cdots,
\mathrm{ip}^{(1,X)@}_{s_{k-1}}\left(
\mathrm{CH}^{(1)}_{s_{k-1}}\left(\mathfrak{P}_{k-1}\right)\right),
\right.\right.
\\&\qquad\qquad\qquad\qquad
\mathrm{ip}^{(1,X)@}_{t}\left(
T^{(1)'}(M)\right),
\mathrm{ip}^{(1,X)@}_{s_{k+1}}\left(
\mathrm{CH}^{(1)}_{s_{k+1}}\left(\mathfrak{P}_{k+1}\right)\right),
\\&\qquad\qquad\qquad\qquad\qquad\qquad\qquad\quad
\left.\left.
\cdots,
\mathrm{ip}^{(1,X)@}_{s_{\bb{\mathbf{s}}-1}}\left(
\mathrm{CH}^{(1)}_{s_{\bb{\mathbf{s}}-1}}\left(\mathfrak{P}_{\bb{\mathbf{s}}-1}\right)\right)
\right)
\right)
\tag{2}
\\&=
\sigma^{\mathbf{T}_{\Sigma}(X)}\left(
\mathrm{sc}^{(0,1)}_{s_{0}}\left(
\mathrm{ip}^{(1,X)@}_{s_{0}}\left(
\mathrm{CH}^{(1)}_{s_{0}}\left(\mathfrak{P}_{0}\right)\right)\right),
\cdots,
\right.
\\&\qquad\qquad
\mathrm{sc}^{(0,1)}_{s_{k-1}}\left(
\mathrm{ip}^{(1,X)@}_{s_{k-1}}\left(
\mathrm{CH}^{(1)}_{s_{k-1}}\left(\mathfrak{P}_{k-1}\right)\right)\right),
\mathrm{sc}^{(0,1)}_{t}\left(
\mathrm{ip}^{(1,X)@}_{t}\left(
T^{(1)'}(M)\right)\right),
\\&\qquad\qquad\qquad\qquad
\mathrm{sc}^{(0,1)}_{s_{k+1}}\left(
\mathrm{ip}^{(1,X)@}_{s_{k+1}}\left(
\mathrm{CH}^{(1)}_{s_{k+1}}\left(\mathfrak{P}_{k+1}\right)\right)\right),
\\&\qquad\qquad\qquad\qquad\qquad\qquad
\left.
\cdots,
\mathrm{sc}^{(0,1)}_{s_{\bb{\mathbf{s}}-1}}\left(
\mathrm{ip}^{(1,X)@}_{s_{\bb{\mathbf{s}}-1}}\left(
\mathrm{CH}^{(1)}_{s_{\bb{\mathbf{s}}-1}}\left(\mathfrak{P}_{\bb{\mathbf{s}}-1}\right)\right)
\right)
\right)
\tag{3}
\\&=
\sigma^{\mathbf{T}_{\Sigma}(X)}\left(
\mathrm{sc}^{(0,1)}_{s_{0}}\left(
\mathrm{ip}^{(1,X)@}_{s_{0}}\left(
\mathrm{CH}^{(1)}_{s_{0}}\left(\mathfrak{P}_{0}\right)\right)\right),
\cdots,
\right.
\\&\qquad\qquad
\mathrm{sc}^{(0,1)}_{s_{k-1}}\left(
\mathrm{ip}^{(1,X)@}_{s_{k-1}}\left(
\mathrm{CH}^{(1)}_{s_{k-1}}\left(\mathfrak{P}_{k-1}\right)\right)\right),
\mathrm{sc}^{(0,1)}_{t}\left(
\mathrm{ip}^{(1,X)@}_{t}\left(
T^{(1)'}(N)\right)\right),
\\&\qquad\qquad\qquad\qquad
\mathrm{sc}^{(0,1)}_{s_{k+1}}\left(
\mathrm{ip}^{(1,X)@}_{s_{k+1}}\left(
\mathrm{CH}^{(1)}_{s_{k+1}}\left(\mathfrak{P}_{k+1}\right)\right)\right),
\\&\qquad\qquad\qquad\qquad\qquad\qquad
\left.
\cdots,
\mathrm{sc}^{(0,1)}_{s_{\bb{\mathbf{s}}-1}}\left(
\mathrm{ip}^{(1,X)@}_{s_{\bb{\mathbf{s}}-1}}\left(
\mathrm{CH}^{(1)}_{s_{\bb{\mathbf{s}}-1}}\left(\mathfrak{P}_{\bb{\mathbf{s}}-1}\right)\right)
\right)
\right)
\tag{4}
\\&=
\mathrm{sc}^{(0,1)}_{s}\left(
\sigma^{\mathbf{Pth}_{\boldsymbol{\mathcal{A}}}}\left(
\mathrm{ip}^{(1,X)@}_{s_{0}}\left(
\mathrm{CH}^{(1)}_{s_{0}}\left(\mathfrak{P}_{0}\right)\right),
\cdots,
\mathrm{ip}^{(1,X)@}_{s_{k-1}}\left(
\mathrm{CH}^{(1)}_{s_{k-1}}\left(\mathfrak{P}_{k-1}\right)\right),
\right.\right.
\\&\qquad\qquad\qquad\qquad
\mathrm{ip}^{(1,X)@}_{t}\left(
T^{(1)'}(N)\right),
\mathrm{ip}^{(1,X)@}_{s_{k+1}}\left(
\mathrm{CH}^{(1)}_{s_{k+1}}\left(\mathfrak{P}_{k+1}\right)\right),
\\&\qquad\qquad\qquad\qquad\qquad\qquad\qquad\quad
\left.\left.
\cdots,
\mathrm{ip}^{(1,X)@}_{s_{\bb{\mathbf{s}}-1}}\left(
\mathrm{CH}^{(1)}_{s_{\bb{\mathbf{s}}-1}}\left(\mathfrak{P}_{\bb{\mathbf{s}}-1}\right)\right)
\right)
\right)
\tag{5}
\\&=
\mathrm{sc}^{(0,1)}_{s}\left(
\mathrm{ip}^{(1,X)@}_{s}\left(
\sigma^{\mathbf{T}_{\Sigma^{\boldsymbol{\mathcal{A}}}}(X)}\left(
\mathrm{CH}^{(1)}_{s_{0}}\left(\mathfrak{P}_{0}\right),
\cdots,
\mathrm{CH}^{(1)}_{s_{k-1}}\left(\mathfrak{P}_{k-1}\right),
\right.\right.\right.
\\&\qquad\qquad\qquad\qquad
\left.\left.\left.
T^{(1)'}(N),
\mathrm{CH}^{(1)}_{s_{k+1}}\left(\mathfrak{P}_{k+1}\right),
\cdots,
\mathrm{CH}^{(1)}_{s_{\bb{\mathbf{s}}-1}}\left(\mathfrak{P}_{\bb{\mathbf{s}}-1}\right)
\right)
\right)\right)
\tag{6}
\\&=\mathrm{sc}^{(0,1)}_{s}\left(\mathrm{ip}^{(1,X)@}_{s}\left(
T^{(1)}(N)\right)\right).
\tag{7}
\end{align*}
\end{flushleft}

In the just stated chain of equalities, the first equality unravels the description of the  first-order elementary translation $T^{(1)}$; the second equality follows from the fact that $\mathrm{ip}^{(1,X)@}$ is a $\Sigma^{\boldsymbol{\mathcal{A}}}$-homomorphism, according to Definition~\ref{DIp}; the third equality follows from the fact that $\mathrm{sc}^{(0,1)}$ is a $\Sigma$-homomorphism, according to Proposition~\ref{PHom}; the fourth equality follows from the fact that, by induction, $
\mathrm{sc}^{(0,1)}_{t}(
\mathrm{ip}^{(1,X)@}_{t}(
T^{(1)'}(M)))
=
\mathrm{sc}^{(0,1)}_{t}(
\mathrm{ip}^{(1,X)@}_{t}(
T^{(1)'}(N)))$;
the fifth equality follows from the fact that $\mathrm{sc}^{(0,1)}$ is a $\Sigma$-homomorphism, according to Proposition~\ref{PHom}; the sixth equality follows from the fact that $\mathrm{ip}^{(1,X)@}$ is a $\Sigma^{\boldsymbol{\mathcal{A}}}$-homomorphism, according to Definition~\ref{DIp}; finally, the last equality recovers  the description of the first-order translation $T^{(1)}$.

The case of the $(0,1)$-target follows by a similar argument.

Now, regarding the other options for the first-order  translation $T^{(1)}$, it could be the case that $s=t$ and there is a path $\mathfrak{P}\in\mathrm{Pth}_{\boldsymbol{\mathcal{A}},s}$ and a first-order translation $T^{(1)'}$ in $\mathrm{Tl}_{t}(\mathrm{T}_{\Sigma^{\boldsymbol{\mathcal{A}}}}(X))_{s}$ of height $m$  such that  such that
\[
T^{(1)}=T^{(1)'}\circ^{0\mathbf{T}_{\Sigma^{\boldsymbol{\mathcal{A}}}}(X)}_{s}\mathrm{CH}^{(1)}_{s}\left(\mathfrak{P}\right).
\]

For this case, the following chain of equivalences hold
\begin{flushleft}
$T^{(1)}(M)$ is a path term
\allowdisplaybreaks
\begin{align*}
\Leftrightarrow\quad & T^{(1)'}(M)\circ^{0\mathbf{T}_{\Sigma^{\boldsymbol{\mathcal{A}}}}(X)}_{s}\mathrm{CH}^{(1)}_{s}\left(\mathfrak{P}\right)\mbox{ is a path term}
\tag{1}
\\
\Leftrightarrow\quad &
\mathrm{ip}^{(1,X)@}_{s}\left(
T^{(1)'}(M)\circ^{0\mathbf{T}_{\Sigma^{\boldsymbol{\mathcal{A}}}}(X)}_{s}\mathrm{CH}^{(1)}_{s}\left(\mathfrak{P}\right)
\right)
\mbox{ is a path}
\tag{2}
\\
\Leftrightarrow\quad &
\mathrm{ip}^{(1,X)@}_{s}\left(
T^{(1)'}(M)\right)\circ^{0\mathbf{Pth}_{\boldsymbol{\mathcal{A}}}}_{s}
\mathrm{ip}^{(1,X)@}_{s}\left(
\mathrm{CH}^{(1)}_{s}\left(\mathfrak{P}\right)
\right)
\mbox{ is a path}
\tag{3}
\\
\Leftrightarrow\quad &
\begin{cases}
\mathrm{ip}^{(1,X)@}_{s}\left(
T^{(1)'}(M)\right)\mbox{ is a path}
\\
\mathrm{sc}^{(0,1)}_{s}\left(
\mathrm{ip}^{(1,X)@}_{s}\left(
T^{(1)'}(M)\right)\right)
=
\mathrm{tg}^{(0,1)}_{s}\left(
\mathrm{ip}^{(1,X)@}_{s}\left(
\mathrm{CH}^{(1)}_{s}\left(\mathfrak{P}\right)
\right)\right)
\end{cases}
\tag{4}
\\
\Leftrightarrow\quad &
\begin{cases}
T^{(1)'}(M)\mbox{ is a path term}
\\
\mathrm{sc}^{(0,1)}_{s}\left(
\mathrm{ip}^{(1,X)@}_{s}\left(
T^{(1)'}(M)\right)\right)
=
\mathrm{tg}^{(0,1)}_{s}\left(
\mathrm{ip}^{(1,X)@}_{s}\left(
\mathrm{CH}^{(1)}_{s}\left(\mathfrak{P}\right)
\right)\right)
\end{cases}
\tag{5}
\\
\Leftrightarrow\quad &
\begin{cases}
T^{(1)'}(N)\mbox{ is a path term}
\\
\mathrm{sc}^{(0,1)}_{s}\left(
\mathrm{ip}^{(1,X)@}_{s}\left(
T^{(1)'}(N)\right)\right)
=
\mathrm{tg}^{(0,1)}_{s}\left(
\mathrm{ip}^{(1,X)@}_{s}\left(
\mathrm{CH}^{(1)}_{s}\left(\mathfrak{P}\right)
\right)\right)
\end{cases}
\tag{6}
\\
\Leftrightarrow\quad &
\begin{cases}
\mathrm{ip}^{(1,X)@}_{s}\left(
T^{(1)'}(N)\right)\mbox{ is a path}
\\
\mathrm{sc}^{(0,1)}_{s}\left(
\mathrm{ip}^{(1,X)@}_{s}\left(
T^{(1)'}(N)\right)\right)
=
\mathrm{tg}^{(0,1)}_{s}\left(
\mathrm{ip}^{(1,X)@}_{s}\left(
\mathrm{CH}^{(1)}_{s}\left(\mathfrak{P}\right)
\right)\right)
\end{cases}
\tag{7}
\\
\Leftrightarrow\quad &
\mathrm{ip}^{(1,X)@}_{s}\left(
T^{(1)'}(N)\right)\circ^{0\mathbf{Pth}_{\boldsymbol{\mathcal{A}}}}_{s}
\mathrm{ip}^{(1,X)@}_{s}\left(
\mathrm{CH}^{(1)}_{s}\left(\mathfrak{P}\right)
\right)
\mbox{ is a path}
\tag{8}
\\
\Leftrightarrow\quad &
\mathrm{ip}^{(1,X)@}_{s}\left(
T^{(1)'}(N)\circ^{0\mathbf{T}_{\Sigma^{\boldsymbol{\mathcal{A}}}}(X)}_{s}\mathrm{CH}^{(1)}_{s}\left(\mathfrak{P}\right)
\right)
\mbox{ is a path}
\tag{9}
\\
\Leftrightarrow\quad & T^{(1)'}(N)\circ^{0\mathbf{T}_{\Sigma^{\boldsymbol{\mathcal{A}}}}(X)}_{s}\mathrm{CH}^{(1)}_{s}\left(\mathfrak{P}\right)\mbox{ is a path term}
\tag{10}
\\
\Leftrightarrow\quad &
T^{(1)}(N) \mbox{ is a path term}.
\tag{11}
\end{align*}
\end{flushleft}

In the just stated chain of equivalences, the first equivalence follows from the description of the first-order  translation $T^{(1)}$; the second equivalence follows from Proposition~\ref{PPT}; the third equivalence follows from the fact that $\mathrm{ip}^{(1,X)@}$ is a $\Sigma^{\boldsymbol{\mathcal{A}}}$-homomorphism, according to Definition~\ref{DIp}; the fourth equivalence follows from the definition of the $0$-composition operation in $\mathbf{Pth}_{\boldsymbol{\mathcal{A}}}$ according to Proposition~\ref{PPthCatAlg}; 
the fifth equivalence follows from Proposition~\ref{PPT};
the  sixth equivalence follows by induction. Moreover, the equality $
\mathrm{sc}^{(0,1)}_{t}(
\mathrm{ip}^{(1,X)@}_{t}(
T^{(1)'}(M)))
=
\mathrm{sc}^{(0,1)}_{t}(
\mathrm{ip}^{(1,X)@}_{t}(
T^{(1)'}(N)))$ also holds; the seventh equivalence follows from Proposition~\ref{PPT}; the eighth equivalence follows from the definition of the $0$-composition operation in $\mathbf{Pth}_{\boldsymbol{\mathcal{A}}}$ according to Proposition~\ref{PPthCatAlg}; the ninth equivalence follows from the fact that $\mathrm{ip}^{(1,X)@}$ is a $\Sigma^{\boldsymbol{\mathcal{A}}}$-homomorphism, according to Definition~\ref{DIp}; the tenth equivalence follows from Proposition~\ref{PPT}; finally, the last equivalence follows from the description of the first-order translation $T^{(1)}$.

Assume that either $T^{(1)}(M)$ or $T^{(1)}(N)$ is a path term.
Regarding the $(0,1)$-source, the following chain of equalities holds
\begin{flushleft}
$\mathrm{sc}^{(0,1)}_{s}\left(\mathrm{ip}^{(1,X)@}_{s}\left(
T^{(1)}(M)\right)\right)$
\allowdisplaybreaks
\begin{align*}
\qquad&=\mathrm{sc}^{(0,1)}_{s}\left(\mathrm{ip}^{(1,X)@}_{s}\left(
T^{(1)'}(M)\circ^{0\mathbf{T}_{\Sigma^{\boldsymbol{\mathcal{A}}}}(X)}_{s}\mathrm{CH}^{(1)}_{s}\left(\mathfrak{P}\right)\right)\right)
\tag{1}
\\&=
\mathrm{sc}^{(0,1)}_{s}\left(
\mathrm{ip}^{(1,X)@}_{s}\left(
T^{(1)'}(M)\right)
\circ^{0\mathbf{Pth}_{\boldsymbol{\mathcal{A}}}}_{s}
\mathrm{ip}^{(1,X)@}_{s}\left(
\mathrm{CH}^{(1)}_{s}\left(
\mathfrak{P}\right)
\right)\right)
\tag{2}
\\&=
\mathrm{sc}^{(0,1)}_{s}\left(
\mathrm{ip}^{(1,X)@}_{s}\left(
\mathrm{CH}^{(1)}_{s}\left(
\mathfrak{P}\right)
\right)\right)
\tag{3}
\\&=
\mathrm{sc}^{(0,1)}_{s}\left(
\mathrm{ip}^{(1,X)@}_{s}\left(
T^{(1)'}(N)\right)
\circ^{0\mathbf{Pth}_{\boldsymbol{\mathcal{A}}}}_{s}
\mathrm{ip}^{(1,X)@}_{s}\left(
\mathrm{CH}^{(1)}_{s}\left(
\mathfrak{P}\right)
\right)\right)
\tag{4}
\\&=\mathrm{sc}^{(0,1)}_{s}\left(\mathrm{ip}^{(1,X)@}_{s}\left(
T^{(1)'}(N)\circ^{0\mathbf{T}_{\Sigma^{\boldsymbol{\mathcal{A}}}}(X)}_{s}\mathrm{CH}^{(1)}_{s}\left(\mathfrak{P}\right)\right)\right)
\tag{5}
\\&=
\mathrm{sc}^{(0,1)}_{s}\left(\mathrm{ip}^{(1,X)@}_{s}\left(
T^{(1)}(N)\right)\right).
\tag{6}
\end{align*}
\end{flushleft}

In the just stated chain of equalities, the first equality follows from the description of the first-order  translation $T^{(1)}$; the second equality follows from the fact that $\mathrm{ip}^{(1,X)@}$ is a $\Sigma^{\boldsymbol{\mathcal{A}}}$-homomorphism, according to Definition~\ref{DIp}; the third equality follows from the definition of the $0$-composition operation in $\mathbf{Pth}_{\boldsymbol{\mathcal{A}}}$ according to Proposition~\ref{PPthCatAlg}; the fourth equality follows from the definition of the $0$-composition operation in $\mathbf{Pth}_{\boldsymbol{\mathcal{A}}}$ according to Proposition~\ref{PPthCatAlg}; the fifth equality follows from the fact that $\mathrm{ip}^{(1,X)@}$ is a $\Sigma^{\boldsymbol{\mathcal{A}}}$-homomorphism, according to Definition~\ref{DIp}; finally, the last equivalence follows from the description of the first-order elementary translation $T^{(1)}$.

Now, regarding the $(0,1)$-target, the following chain of equalities holds
\begin{flushleft}
$\mathrm{tg}^{(0,1)}_{s}\left(\mathrm{ip}^{(1,X)@}_{s}\left(
T^{(1)}(M)\right)\right)$
\allowdisplaybreaks
\begin{align*}
\qquad&=\mathrm{tg}^{(0,1)}_{s}\left(\mathrm{ip}^{(1,X)@}_{s}\left(
T^{(1)'}(M)\circ^{0\mathbf{T}_{\Sigma^{\boldsymbol{\mathcal{A}}}}(X)}_{s}\mathrm{CH}^{(1)}_{s}\left(\mathfrak{P}\right)\right)\right)
\tag{1}
\\&=
\mathrm{tg}^{(0,1)}_{s}\left(
\mathrm{ip}^{(1,X)@}_{s}\left(
T^{(1)'}(M)\right)
\circ^{0\mathbf{Pth}_{\boldsymbol{\mathcal{A}}}}_{s}
\mathrm{ip}^{(1,X)@}_{s}\left(
\mathrm{CH}^{(1)}_{s}\left(
\mathfrak{P}\right)
\right)\right)
\tag{2}
\\&=
\mathrm{tg}^{(0,1)}_{s}\left(
\mathrm{ip}^{(1,X)@}_{s}\left(
T^{(1)'}(M)\right)
\right)
\tag{3}
\\&=
\mathrm{tg}^{(0,1)}_{s}\left(
\mathrm{ip}^{(1,X)@}_{s}\left(
T^{(1)'}(N)\right)
\right)
\tag{4}
\\&=
\mathrm{tg}^{(0,1)}_{s}\left(
\mathrm{ip}^{(1,X)@}_{s}\left(
T^{(1)'}(N)\right)
\circ^{0\mathbf{Pth}_{\boldsymbol{\mathcal{A}}}}_{s}
\mathrm{ip}^{(1,X)@}_{s}\left(
\mathrm{CH}^{(1)}_{s}\left(
\mathfrak{P}\right)
\right)\right)
\tag{5}
\\&=\mathrm{tg}^{(0,1)}_{s}\left(\mathrm{ip}^{(1,X)@}_{s}\left(
T^{(1)'}(N)\circ^{0\mathbf{T}_{\Sigma^{\boldsymbol{\mathcal{A}}}}(X)}_{s}\mathrm{CH}^{(1)}_{s}\left(\mathfrak{P}\right)\right)\right)
\tag{6}
\\&=
\mathrm{tg}^{(0,1)}_{s}\left(\mathrm{ip}^{(1,X)@}_{s}\left(
T^{(1)}(N)\right)\right).
\tag{7}
\end{align*}
\end{flushleft}

In the just stated chain of equalities, the first equality follows from the description of the first-order  translation $T^{(1)}$; the second equality follows from the fact that $\mathrm{ip}^{(1,X)@}$ is a $\Sigma^{\boldsymbol{\mathcal{A}}}$-homomorphism, according to Definition~\ref{DIp}; the third equality follows from the definition of the $0$-composition operation in $\mathbf{Pth}_{\boldsymbol{\mathcal{A}}}$ according to Proposition~\ref{PPthCatAlg}; the fourth equality follows by induction. Note that, according to the above proof,  since $T^{(1)}(M)$ is a path term, then $T^{(1)'}(M)$ is a path term as well, which implies, by induction, that  $T^{(1)'}(N)$ is a path term. The implication also holds in the other direction. In particular,  $
\mathrm{tg}^{(0,1)}_{t}(
\mathrm{ip}^{(1,X)@}_{t}(
T^{(1)'}(M)))
=
\mathrm{tg}^{(0,1)}_{t}(
\mathrm{ip}^{(1,X)@}_{t}(
T^{(1)'}(N)))$; the fifth equality follows from the definition of the $0$-composition operation in $\mathbf{Pth}_{\boldsymbol{\mathcal{A}}}$ according to Proposition~\ref{PPthCatAlg}; the sixth equality follows from the fact that $\mathrm{ip}^{(1,X)@}$ is a $\Sigma^{\boldsymbol{\mathcal{A}}}$-homomorphism, according to Definition~\ref{DIp}; finally, the last equivalence follows from the description of the first-order  translation $T^{(1)}$.

The remaining case, i.e., the case for which  $s=t$ and there is a path $\mathfrak{P}\in\mathrm{Pth}_{\boldsymbol{\mathcal{A}},s}$ and a first-order translation $T^{(1)'}$ in $\mathrm{Tl}_{t}(\mathrm{T}_{\Sigma^{\boldsymbol{\mathcal{A}}}}(X))_{s}$ of height $m$  such that 
\[
T^{(1)}=\mathrm{CH}^{(1)}_{s}\left(\mathfrak{P}\right)\circ^{0\mathbf{T}_{\Sigma^{\boldsymbol{\mathcal{A}}}}(X)}_{s} T^{(1)'},
\]
follows by using an argument similar to the one used in the proof presented above.

This complete the proof.
\end{proof}

\begin{restatable}{proposition}{PUTransWD}
\label{PUTransWD} Let $s$ and $t$ be sorts in $S$, $T^{(1)}$ a first-order translation in 
$\mathrm{Tl}_{t}(\mathrm{T}_{\Sigma^{\boldsymbol{\mathcal{A}}}}(X))_{s}$ and $P$ and $P'$ path terms in
$\mathrm{PT}_{\boldsymbol{\mathcal{A}}, t}$ such that $(P,P')\in \Theta^{[1]}_{t}$. If either $T^{(1)}(P)$ or $T^{(1)}(P')$ is a path term in $\mathrm{PT}_{\boldsymbol{\mathcal{A}}, s}$, then 
\[
\left(
T^{(1)}(P),T^{(1)}(P')
\right)
\in\Theta^{[1]}_{s}.
\]
\end{restatable}
\begin{proof}
It follows from the fact that $\Theta^{[1]}$ is a $\Sigma^{\boldsymbol{\mathcal{A}}}$-congruence on $\mathbf{T}_{\Sigma^{\boldsymbol{\mathcal{A}}}}$, by Definition~\ref{DThetaCong}, and by Proposition~\ref{PTransCong}.
\end{proof}


\part{Second-order rewriting systems}
\chapter{
\texorpdfstring
{Second-order paths on path terms}
{Second-order paths}
}

In this chapter we introduce the notions of second-order rewrite rule and second-order rewrite system $\boldsymbol{\mathcal{A}}^{(2)}$. Following the result of the Chapter~16, the second-order rewrite rules require a condition on the $(0,1)$-sources and $(0,1)$-targets so that the first-order translations are well-defined. This allow us to introduce the concept of second-order path between path term classes. Then we consider the many-sorted set of second-order paths over $\boldsymbol{\mathcal{A}}^{(2)}$, denoted by $\mathrm{Pth}_{\boldsymbol{\mathcal{A}}^{(2)}}$. We further define the notions of $([1],2)$-source, $([1],2)$-target and $(2,[1])$-identity second-order path, denoted by $\mathrm{sc}^{([1],2)}$, $\mathrm{tg}^{([1],2)}$ and $\mathrm{ip}^{(2,[1])\sharp}$, respectively. We then start the structural study of these objects by characterising the $(2,[1])$-identity second-order paths. We also present the $1$-composition of second-order paths and prove that this partial operation is well-defined. In the following we introduce the notion of subpath and we relate the concepts of $1$-composition and subpaths. We then go on to introduce the notion of the second-order echelon of a second-order rewriting rule and the concepts of  echelonless second-order path and head-constant second-order paths. Contrary to what happens in the first part, echelonless second-order paths are not necessarily head-constant. This requires to introduce the notion of coherent head-constant echelonless second-order path. We prove that one-step echelonless second-order paths are coherent and head-constant. Also, head-constant echelonless second-order paths associated to an operation symbol in the original signature $\Sigma$ are coherent. We conclude this chapter by introducing the extraction algorithm for coherent head-constant echelonless second-order paths.


\begin{restatable}{definition}{DDRewSys}
\label{DDRewSys}
A \emph{ second-order many-sorted rewriting system}\index{rewriting system!second-order, $\boldsymbol{\mathcal{A}}^{(2)}$} (or, simply, a \emph{second-order rewriting system}) is an ordered quintuple
$
(S,\Sigma,X,\mathcal{A},\mathcal{A}^{(2)}),
$
often abbreviated to $\boldsymbol{\mathcal{A}}^{(2)}$, where $(S,\Sigma,X,\mathcal{A})$ is a many-sorted rewriting system, see Definition~\ref{DRewSys}, and, for the many-sorted signature $\mathbf{\Sigma}^{\boldsymbol{\mathcal{A}}} = (S,\Sigma^{\boldsymbol{\mathcal{A}}})$, $\mathcal{A}^{(2)}$ a subset of 
\begin{multline*}
\mathrm{Rwr}(\mathbf{\Sigma}^{\boldsymbol{\mathcal{A}}}, X)
=
\left(
\left\lbrace
([M]_{s}, [N]_{s}) \in [\mathrm{PT}_{\boldsymbol{\mathcal{A}}}]^{2}_{s} 
\middle|
\right.\right.
\\
\left.\left.
\left(
\mathrm{ip}^{(1,X)@}_{s}(M),
\mathrm{ip}^{(1,X)@}_{s}(N)
\right)
\in \mathrm{Ker}\left(\mathrm{sc}^{(0,1)}\right)_{s}
\cap\mathrm{Ker}\left(\mathrm{tg}^{(0,1)}\right)_{s}
\right\rbrace
\right)_{s\in S},
\end{multline*}
the $S$-sorted set of the \emph{second-order rewrite rules with variables in $X$}.

\index{rewrite rule!second-order, $\mathfrak{p}^{(2)}$}
For $s\in S$,  we will call the elements of $\mathrm{Rwr}(\mathbf{\Sigma}^{\boldsymbol{\mathcal{A}}}, X)_{s}$ \emph{second-order rewrite rules of type $(X,s)$} and we will denote them with lowercase Euler fraktur letters with the superscript $(2)$, indicating the order, with or without subscripts, e.g., $\mathfrak{p}^{(2)}$, $\mathfrak{p}^{(2)}_{i}$, $\mathfrak{q}^{(2)}$, $\mathfrak{q}^{(2)}_{i}$, etc). 

We will say that $\boldsymbol{\mathcal{A}}^{(2)}$ is a \emph{finite} second-order rewriting system if $\mathcal{A}^{(2)}$ is finite, i.e., if $\mathrm{supp}_{S}(\mathcal{A}^{(2)})$ is finite and, for every $s\in \mathrm{supp}_{S}(\mathcal{A}^{(2)})$, $\mathcal{A}^{(2)}_{s}$ is finite.
\end{restatable}

\begin{assumption}
From now on $\boldsymbol{\mathcal{A}}^{(2)}$ stands for a second-order many-sorted rewriting system, fixed once and for all.
\end{assumption}

We next define the notion of second-order path in $\boldsymbol{\mathcal{A}}^{(2)}$ from a path term class to another.

\begin{restatable}{definition}{DDPth}
\label{DDPth} Let $s$ be a sort in $S$, $[P]_{s},[Q]_{s}$ path term classes in $[\mathrm{PT}_{\boldsymbol{\mathcal{A}}}]_{s}$ and $\mathbf{c}$ a word in $S^{\star}$. Then a \emph{second-order $\mathbf{c}$-path in $\boldsymbol{\mathcal{A}}^{(2)}$ from $[P]_{s}$ to $[Q]_{s}$} is an ordered triple
\index{path!second-order!$\mathfrak{P}^{(2)}$}
$$
\mathfrak{P}^{(2)}
=
\left(
(
[P_{i}]_{s}
)_{i\in\bb{\mathbf{c}}+1}, 
(
\mathfrak{p}^{(2)}_{i}
)_{i\in\bb{\mathbf{c}}},
(
T^{(1)}_{i}
)_{i\in\bb{\mathbf{c}}}
\right)
$$
in $
[\mathrm{PT}_{\boldsymbol{\mathcal{A}}}]^{\bb{\mathbf{c}}+1}_{s}
\times
\mathcal{A}^{(2)}_{\mathbf{c}}
\times
\mathrm{Tl}_{\mathbf{c}}(\mathbf{T}_{\Sigma^{\boldsymbol{\mathcal{A}}}}(X))_{s}$, where $\mathcal{A}^{(2)}_{\mathbf{c}} = \prod_{i\in\bb{\mathbf{c}}}\mathcal{A}^{(2)}_{c_{i}}$ and $\mathrm{Tl}_{\mathbf{c}}(\mathbf{T}_{\Sigma^{\boldsymbol{\mathcal{A}}}}(X))_{s} = \prod_{i\in\bb{\mathbf{c}}}\mathrm{Tl}_{c_{i}}(\mathbf{T}_{\Sigma^{\boldsymbol{\mathcal{A}}}}(X))_{s}$, with $\mathrm{Tl}_{c_{i}}(\mathbf{T}_{\Sigma^{\boldsymbol{\mathcal{A}}}}(X))_{s}$ the set of the first-order $c_{i}$-translations of sort $c$ for $\mathbf{T}_{\Sigma^{\boldsymbol{\mathcal{A}}}}(X)$, see Definition~\ref{DUTrans}, 
such that
\begin{itemize}
\item[(1)] $[P_{0}]_{s}=[P]_{s}$,
\item[(2)] $[P_{\bb{\mathbf{c}}}]_{s}=[Q]_{s}$, and
\item[(3)] for every $i\in\bb{\mathbf{c}}$, if $\mathfrak{p}^{(2)}_{i}=
(
[M_{i}]_{c_{i}}, 
[N_{i}]_{c_{i}}
),
$  
then
\begin{itemize}
\item[(i)] $T^{(1)}_{i}\left(
M_{i}
\right)\in [P_{i}]_{s}$ and
\item[(ii)] $T^{(1)}_{i}\left(
N_{i}
\right)\in [P_{i+1}]_{s}$.
\end{itemize}
\end{itemize}
That is, at each step $i$, we consider a second-order rewrite rule $\mathfrak{p}^{(2)}_{i}$ and a first-order $c_{i}$-translation of sort $s$ $T^{(1)}_{i}$ for $\mathbf{T}_{\Sigma^{\boldsymbol{\mathcal{A}}}}(X)$, see Definition~\ref{DUTrans}, and we require that (i) the value of $T^{(1)}_{i}$ at $M_{i}$ is in $[P_{i}]_{s}$ and (ii) the value of $T^{(1)}_{i}$ at $N_{i}$ is in $[P_{i+1}]_{s}$. In this regard, let us recall that, by Lemma~\ref{LUTransWD}, the condition on the $(0,1)$-source and $(0,1)$-target we have imposed in Definition~\ref{DDRewSys} guarantees that the values of the first-order translations will always be path terms. Moreover, according to Proposition~\ref{PUTransWD}, ultimately, it does not matter which representative term we use.

These second-order paths will be variously depicted as 
$\mathfrak{P}^{(2)}\colon [P]_{s}{\implies}[Q]_{s}$, 
$\mathfrak{P}^{(2)}\colon [P_{0}]_{s}{\implies}[P_{\bb{\mathbf{c}}}]_{s}$,
or
$$
\xymatrixrowsep{15pt}
\xymatrix@C=80pt{
\mathfrak{P}^{(2)}: [P_{0}]_{s}
\ar@{=>}[r]^-{\text{\Small{($\mathfrak{p}^{(2)}_{0}$, $T^{(1)}_{0}$)}}}
&
\quad
[P_{1}]_{s}\quad
\ar@{=>}[r]^-{\text{\Small{($\mathfrak{p}^{(2)}_{1}$, $T^{(1)}_{1}$)}}}
&
{}
\hdots\qquad
{}
\\
\qquad\hdots
{}
\ar@{=>}[r]^-{\text{\Small{($\mathfrak{p}^{(2)}_{\bb{\mathbf{c}}-2}$, $T^{(1)}_{\bb{\mathbf{c}}-2}$)}}}
&
[P_{\bb{\mathbf{c}}-1}]_{s}
\ar@{=>}[r]^-{\text{\Small{($\mathfrak{p}^{(2)}_{\bb{\mathbf{c}}-1}$, $T^{(1)}_{\bb{\mathbf{c}}-1}$)}}}
&
[P_{\bb{\mathbf{c}}}]_{s}.
}
$$
(the following variant:
$$
\xymatrixrowsep{15pt}
\xymatrix@C=80pt{
\mathfrak{P}^{(2)}: [P_{0}]_{s}
\ar@{=>}[r]^-{\text{\Small{($\mathfrak{p}^{(2)}$, $T^{(1)})_{0}$}}}
&
\quad
[P_{1}]_{s}
\quad
\ar@{=>}[r]^-{\text{\Small{($\mathfrak{p}^{(2)}$, $T^{(1)})_{1}$}}}
&
{}
\hdots\qquad
{}
\\
\qquad\qquad\hdots
{}
\ar@{=>}[r]^-{\text{\Small{($\mathfrak{p}^{(2)}$, $T^{(1)})_{\bb{\mathbf{c}}-2}$}}}
&
[P_{\bb{\mathbf{c}}-1}]_{s}
\ar@{=>}[r]^-{\text{\Small{($\mathfrak{p}^{(2)}$, $T^{(1)})_{\bb{\mathbf{c}}-1}$}}}
&
[P_{\bb{\mathbf{c}}}]_{s}
}
$$
will be occasionally used). 

For every $i\in \bb{\mathbf{c}}$, we will say that $[P_{i+1}]_{s}$ is
$(\mathfrak{p}^{(2)}_{i}, T^{(1)}_{i})$-\emph{directly derivable} or, when no confusion can arise, \emph{directly derivable}\index{directly derivable} from $[P_{i}]_{s}$. For every $i\in \bb{\mathbf{c}}+1$, the path term class $[P_{i}]_{s}$ will be called a \emph{$1$-constituent}\index{constituent!$1$-constituent} of the second-order $\bb{\mathbf{c}}$-path $\mathfrak{P}^{(2)}$. Let us note that all $1$-constituents of a second-order path $\mathfrak{P}^{(2)}$ have the same sort. The path term class $[P_{0}]_{s}$ will be called the \emph{$([1],2)$-source}\index{source!second-order!$\mathrm{sc}^{([1],2)}$} of the second-order $\mathbf{c}$-path $\mathfrak{P}^{(2)}$, the path term class $[P_{\bb{\mathbf{c}}}]_{s}$ will be called the \emph{$([1],2)$-target}\index{target!second-order!$\mathrm{tg}^{([1],2)}$} of the second-order $\mathbf{c}$-path $\mathfrak{P}^{(2)}$, and we will say that $\mathfrak{P}^{(2)}$ is a \emph{second-order path from} $[P_{0}]_{s}$ \emph{to} $[P_{\bb{\mathbf{c}}}]_{s}$. 

The \emph{length}\index{path!second-order!length} of a second-order $\mathbf{c}$-path $\mathfrak{P}^{(2)}$ in $\boldsymbol{\mathcal{A}}^{(2)}$, denoted by $\bb{\mathfrak{P}^{(2)}}$, is $\bb{\mathbf{c}}$ and we will say that $\mathfrak{P}^{(2)}$ has $\bb{\mathbf{c}}$ \emph{steps}\index{path!second-order!step}.  If $\bb{\mathfrak{P}^{(2)}}=0$, then we will say that $\mathfrak{P}^{(2)}$ is a \emph{$(2,[1])$-identity second-order path}\index{identity!second-order!$\mathrm{ip}^{(2,[1])\sharp}$}. Let us point out that $\mathfrak{P}^{(2)}$ is a $(2,[1])$-identity second-order path if, and only if, there exists a sort $s$ in $S$ and a path term class $[P]_{s}$ in $[\mathrm{PT}_{\boldsymbol{\mathcal{A}}}]_{s}$ such that, for $\lambda\in S^{\star}$, the empty word on $S$, we have that $(([P]_{s}),\lambda,\lambda)$, identified to $([P]_{s},\lambda, \lambda)$, where, by abuse of notation, we have written $(\lambda, \lambda)$ for the unique element of $\mathcal{A}^{(2)}_{\lambda}\times \mathrm{Tl}_{\lambda}(\mathbf{T}_{\Sigma^{\boldsymbol{\mathcal{A}}}}(X))_{s}$, is equal to $\mathfrak{P}^{(2)}$. This path will be called the \emph{$(2,[1])$-identity second-order path on} $[P]_{s}$. 

If $\bb{\mathfrak{P}^{(2)}}=1$, then we will say that $\mathfrak{P}^{(2)}$ is a \emph{one-step second-order path}\index{path!second-order!one-step}. Moreover, in this case, i.e., when $\mathbf{c} = (c)$, for a unique $c\in S$, and identifying, when no confusion can arise, $(c)$ with $c$, we will speak of  second-order $c$-paths, instead of second-order $(c)$-paths.


We will denote by
\begin{enumerate}
\item $\mathrm{Pth}_{\mathbf{c},\boldsymbol{\mathcal{A}}^{(2)},s}(
[P]_{s},[Q]_{s})$ the set of all second-order $\mathbf{c}$-paths in $\boldsymbol{\mathcal{A}}^{(2)}$ from $[P]_{s}$ to $[Q]_{s}$; by
\item $\mathrm{Pth}_{\boldsymbol{\mathcal{A}}^{(2)},s}(P,Q)$ the set $\bigcup_{\mathbf{c}\in S^{\star}}\mathrm{Pth}_{\mathbf{c},\boldsymbol{\mathcal{A}}^{(2)},s}([P]_{s},[Q]_{s})$ and we will call its elements \emph{second-order paths in $\boldsymbol{\mathcal{A}}^{(2)}$ from $[P]_{s}$ to $[Q]_{s}$};
by
\item $\mathrm{Pth}_{\boldsymbol{\mathcal{A}}^{(2)},s}([P]_{s},\cdot)$ the set $\bigcup_{[Q]_{s}\in[\mathrm{PT}_{\boldsymbol{\mathcal{A}}}]_{s}}\mathrm{Pth}_{\boldsymbol{\mathcal{A}}^{(2)},s}([P]_{s},[Q]_{s})$ and we will call its elements \emph{second-order paths in $\boldsymbol{\mathcal{A}}^{(2)}$ from $[P]_{s}$};
by
\item $\mathrm{Pth}_{\boldsymbol{\mathcal{A}}^{(2)},s}(\cdot,[Q]_{s})$ the set $\bigcup_{[P]_{s}\in[\mathrm{PT}_{\boldsymbol{\mathcal{A}}}]_{s}}\mathrm{Pth}_{\boldsymbol{\mathcal{A}}^{(2)},s}([P]_{s},[Q]_{s})$ and we will call its elements \emph{second-order paths in $\boldsymbol{\mathcal{A}}^{(2)}$ to $[Q]_{s}$}; by
\item $\mathrm{Pth}_{\boldsymbol{\mathcal{A}}^{(2)},s}$ the set $\bigcup_{[P]_{s},[Q]_{s}\in[\mathrm{PT}_{\boldsymbol{\mathcal{A}}}]_{s}}
\mathrm{Pth}_{\boldsymbol{\mathcal{A}}^{(2)},s}([P]_{s},[Q]_{s})$; and, finally, by
\item $\mathrm{Pth}_{\boldsymbol{\mathcal{A}}^{(2)}}$ the $S$-sorted set defined, for every $s\in S$, as $\mathrm{Pth}_{\boldsymbol{\mathcal{A}}^{(2)},s}$.
\end{enumerate}
\index{path!second-order!$\mathrm{Pth}_{\boldsymbol{\mathcal{A}}^{(2)}}$}

We will denote by
\begin{enumerate}
\item $\mathrm{ip}^{(2,X)}$ the $S$-sorted mapping from $X$ to $\mathrm{Pth}_{\boldsymbol{\mathcal{A}}^{(2)}}$ that, for every sort $s\in S$, sends a variable $x\in X_{s}$ to $([\eta^{(1,X)}_{s}(x)]_{s},\lambda,\lambda)$ the $(2,[1])$-identity path on $[\eta^{(1,X)}_{s}(x)]_{s}$; by
\item $\mathrm{ech}^{(2,\mathcal{A})}$ the $S$-sorted mapping from $\mathcal{A}$ to $\mathrm{Pth}_{\boldsymbol{\mathcal{A}}^{(2)}}$ that, for every sort $s\in S$, sends a rewrite rule $\mathfrak{p}\in \mathcal{A}_{s}$ to $([\eta^{(1,\mathcal{A})}_{s}(\mathfrak{p})]_{s},\lambda,\lambda)$ the $(2,[1])$-identity path on $[\eta^{(1,\mathcal{A})}_{s}(\mathfrak{p})]_{s}$; by
\item $\mathrm{sc}^{([1],2)}$ the $S$-sorted mapping from $\mathrm{Pth}_{\boldsymbol{\mathcal{A}}^{(2)}}$ to $[\mathrm{PT}_{\boldsymbol{\mathcal{A}}}]$ that sends a second-order path to its $([1],2)$-source; by
\item $\mathrm{tg}^{([1],2)}$ the $S$-sorted mapping from $\mathrm{Pth}_{\boldsymbol{\mathcal{A}}^{(2)}}$ to $[\mathrm{PT}_{\boldsymbol{\mathcal{A}}}]$ that sends a second-order path to its $([1],2)$-target; and by 
\item $\mathrm{ip}^{(2,[1])\sharp}$ the $S$-sorted mapping from $[\mathrm{PT}_{\boldsymbol{\mathcal{A}}}]$ to $\mathrm{Pth}_{\boldsymbol{\mathcal{A}}^{(2)}}$ that sends a path term class to its $(2,[1])$-identity  second-order path. 
\end{enumerate}
These $S$-sorted mappings are depicted in the diagram of Figures~\ref{FDPthX} and~\ref{FDPthA}.

Finally, given a sort $s$ in $S$ and a second-order path $\mathfrak{P}^{(2)}$ in $\mathrm{Pth}_{\boldsymbol{\mathcal{A}}^{(2)},s}$ we will say that $\mathfrak{P}^{(2)}$ is a \emph{second-order $([1],2)$-loop}\index{loop!second-order!$([1],2)$-second-order loop} if $\mathrm{sc}^{([1],2)}_{s}(\mathfrak{P}^{(2)})=\mathrm{tg}^{([1],2)}_{s}(\mathfrak{P}^{(2)})$. Let us note that every $(2,[1])$-identity second-order path is a second-order $([1],2)$-loop.
\end{restatable}

\begin{figure}
\begin{tikzpicture}
[ACliment/.style={-{To [angle'=45, length=5.75pt, width=4pt, round]}},scale=1.1]
\node[] (xoq) at (0,0) [] {$X$};
\node[] (txoq) at (6,0) [] {$[\mathrm{PT}_{\boldsymbol{\mathcal{A}}}]$};
\node[] (p) at (6,-3) [] {$\mathrm{Pth}_{\boldsymbol{\mathcal{A}}^{(2)}}$};
\draw[ACliment]  (xoq) to node [above]
{$\eta^{([1],X)}$} (txoq);
\draw[ACliment, bend right=10]  (xoq) to node [below left] {$\mathrm{ip}^{(2,X)}$} (p);

\node[] (B0) at (6,-1.5)  [] {};
\draw[ACliment]  ($(B0)+(0,1.2)$) to node [above, fill=white] {
$\textstyle \mathrm{ip}^{(2,[1])\sharp}$
} ($(B0)+(0,-1.2)$);
\draw[ACliment, bend right]  ($(B0)+(.3,-1.2)$) to node [ below, fill=white] {
$\textstyle \mathrm{tg}^{([1],2)}$
} ($(B0)+(.3,1.2)$);
\draw[ACliment, bend left]  ($(B0)+(-.3,-1.2)$) to node [below, fill=white] {
$\textstyle \mathrm{sc}^{([1],2)}$
} ($(B0)+(-.3,1.2)$);
\end{tikzpicture}
\caption{Many-sorted mappings relative to $X$ at layers [1] \& 2.}\label{FDPthX}
\end{figure}

\begin{figure}
\begin{tikzpicture}
[ACliment/.style={-{To [angle'=45, length=5.75pt, width=4pt, round]}},scale=1.1]
\node[] (xoq) at (0,0) [] {$\mathcal{A}$};
\node[] (txoq) at (6,0) [] {$[\mathrm{PT}_{\boldsymbol{\mathcal{A}}}]$};
\node[] (p) at (6,-3) [] {$\mathrm{Pth}_{\boldsymbol{\mathcal{A}}^{(2)}}$};
\draw[ACliment]  (xoq) to node [above]
{$\eta^{([1],\mathcal{A})}$} (txoq);
\draw[ACliment, bend right=10]  (xoq) to node [below left] {$\mathrm{ech}^{(2,\mathcal{A})}$} (p);

\node[] (B0) at (6,-1.5)  [] {};
\draw[ACliment]  ($(B0)+(0,1.2)$) to node [above, fill=white] {
$\textstyle \mathrm{ip}^{(2,[1])\sharp}$
} ($(B0)+(0,-1.2)$);
\draw[ACliment, bend right]  ($(B0)+(.3,-1.2)$) to node [ below, fill=white] {
$\textstyle \mathrm{tg}^{([1],2)}$
} ($(B0)+(.3,1.2)$);
\draw[ACliment, bend left]  ($(B0)+(-.3,-1.2)$) to node [below, fill=white] {
$\textstyle \mathrm{sc}^{([1],2)}$
} ($(B0)+(-.3,1.2)$);
\end{tikzpicture}
\caption{Many-sorted mappings relative to $\mathcal{A}$ at layers [1] \& 2.}\label{FDPthA}
\end{figure}

It follows from Definition~\ref{DDPth} that the $(2,[1])$-identity second-order path on a path term class is the only second-order path of length $0$ whose $([1],2)$-source (or $([1],2)$-target) equals this path term class. Consequently, the only interesting properties about $(2,[1])$-identity second-order paths are those having to do with their $([1],2)$-source (or $([1],2)$-target). 

We next characterize the $(2,[1])$-identity second-order paths as fixed points.

\begin{restatable}{proposition}{PDPthId}
\label{PDPthId} Let $s$ be a sort in $S$ and $\mathfrak{P}^{(2)}$ a second-order path in $\mathrm{Pth}_{\boldsymbol{\mathcal{A}}^{(2)},s}$. Then the following statements are equivalent:
\begin{enumerate}
\item $\mathfrak{P}^{(2)}$ is a $(2,[1])$-identity second-order path.
\item $\mathfrak{P}^{(2)}=\mathrm{ip}^{(2,[1])\sharp}_{s}(\mathrm{sc}^{([1],2)}_{s}(\mathfrak{P}^{(2)}))$, i.e., 
$\mathfrak{P}^{(2)}\in \mathrm{Fix}(\mathrm{ip}^{(2,[1])\sharp}_{s}\circ\mathrm{sc}^{([1],2)}_{s})$, the set of fixed points of $\mathrm{ip}^{(2,[1])\sharp}_{s}\circ\mathrm{sc}^{([1],2)}_{s}$.
\item $\mathfrak{P}^{(2)}=\mathrm{ip}^{(2,[1])\sharp}_{s}(\mathrm{tg}^{([1],2)}_{s}(\mathfrak{P}^{(2)}))$, i.e., 
$\mathfrak{P}^{(2)}\in \mathrm{Fix}(\mathrm{ip}^{(2,[1])\sharp}_{s}\circ\mathrm{tg}^{([1],2)}_{s})$, the set of fixed points of $\mathrm{ip}^{(2,[1])\sharp}_{s}\circ\mathrm{tg}^{([1],2)}_{s}$.
\end{enumerate}
\end{restatable}

The following are some of the relationships between the $S$-sorted mappings $\eta^{([1],X)}$, $\mathrm{ech}^{(2,\mathcal{A})}$, $\mathrm{ip}^{(2,X)}$, $\mathrm{sc}^{([1],2)}$, $\mathrm{tg}^{([1],2)}$,  $\mathrm{ip}^{(2,[1])\sharp}$ and $\mathrm{id}^{[\mathrm{PT}_{\boldsymbol{\mathcal{A}}}]}$.

\begin{proposition}\label{PDBasicEq} The following equalities hold
\begin{itemize}
\item[(i)] $\mathrm{sc}^{([1],2)}\circ\mathrm{ip}^{(2,[1])\sharp}=\mathrm{id}^{[\mathrm{PT}_{\boldsymbol{\mathcal{A}}}]};$
\item[(ii)] $\mathrm{tg}^{([1],2)}\circ\mathrm{ip}^{(2,[1])\sharp}=\mathrm{id}^{[\mathrm{PT}_{\boldsymbol{\mathcal{A}}}]}.$
\end{itemize}
The reader is advised to consult the diagram of Figure~\ref{FDPthX}.
\end{proposition}

\begin{proposition}\label{PDBasicEqX} 
The following equalities hold
\begin{itemize}
\item[(i)] $\mathrm{sc}^{([1],2)}\circ\mathrm{ip}^{(2,X)}=\eta^{([1],X)};$
\item[(ii)] $\mathrm{tg}^{([1],2)}\circ\mathrm{ip}^{(2,X)}=\eta^{([1],X)};$
\item[(iii)] $\mathrm{ip}^{(2,[1])\sharp}\circ \eta^{([1],X)}=\mathrm{ip}^{(2,X)}.$
\end{itemize}
The reader is advised to consult the diagram of Figure~\ref{FDPthX}.
\end{proposition}

\begin{proposition}\label{PDBasicEqA} 
The following equalities hold
\begin{itemize}
\item[(i)] $\mathrm{sc}^{([1],2)}\circ\mathrm{ech}^{(2,\mathcal{A})}=\eta^{([1],\mathcal{A})};$
\item[(ii)] $\mathrm{tg}^{([1],2)}\circ\mathrm{ech}^{(2,\mathcal{A})}=\eta^{([1],\mathcal{A})};$
\item[(iii)] $\mathrm{ip}^{(2,[1])\sharp}\circ \eta^{([1],\mathcal{A})}=\mathrm{ech}^{(2,\mathcal{A})}.$
\end{itemize}
The reader is advised to consult the diagram of Figure~\ref{FDPthA}.
\end{proposition}

\begin{remark}
Later on, after defining a suitable structure of $\Sigma^{\boldsymbol{\mathcal{A}}}$-algebra on $\mathrm{Pth}_{\boldsymbol{\mathcal{A}}^{(2)}}$, which will give raise to the $\Sigma^{\boldsymbol{\mathcal{A}}}$-algebra $\mathbf{Pth}^{(1,2)}_{\boldsymbol{\mathcal{A}}^{(2)}}$, we will show that $\mathrm{sc}^{([1],2)}$ and $\mathrm{tg}^{([1],2)}$ are homomorphisms from $\mathbf{Pth}^{(1,2)}_{\boldsymbol{\mathcal{A}}^{(2)}}$ to $[\mathbf{PT}_{\boldsymbol{\mathcal{A}}}]$, and  that $\mathrm{ip}^{(2,[1])\sharp}$ is, in fact, the unique homomorphism from $[\mathbf{PT}_{\boldsymbol{\mathcal{A}}}]$ to $\mathbf{Pth}^{(1,2)}_{\boldsymbol{\mathcal{A}}^{(2)}}$ such that $\mathrm{ip}^{(2,[1])} = \mathrm{ip}^{(2,[1])\sharp}\circ \eta^{([1],X)}$.
\end{remark}

We next define, for every sort $s$ in $S$, the partial operation of $[1]$-composition of second-order paths of sort $s$.

\begin{restatable}{definition}{DDPthComp}
\label{DDPthComp}
Let $s$ be a sort in $S$, $\mathfrak{P}^{(2)}$ and $\mathfrak{Q}^{(2)}$ second-order paths in $\mathrm{Pth}_{\boldsymbol{\mathcal{A}}^{(2)},s}$, where, for a unique word $\mathbf{c}\in S^{\star}$, $\mathfrak{P}^{(2)}$ is a second-order $\mathbf{c}$-path in $\boldsymbol{\mathcal{A}}^{(2)}$ of the form
$$
\mathfrak{P}^{(2)}
=
\left(
([P_{i}]_{s})_{i\in\bb{\mathbf{c}}+1},
(\mathfrak{p}^{(2)}_{i})_{i\in\bb{\mathbf{c}}},
(T^{(1)}_{i})_{i\in\bb{\mathbf{c}}}
\right),
$$
and, for a unique $\mathbf{d}\in S^{\star}$, $\mathfrak{Q}^{(2)}$ is a second-order $\mathbf{d}$-path in $\boldsymbol{\mathcal{A}}^{(2)}$ of the form
$$
\mathfrak{Q}^{(2)}
=
\left(
([Q_{j}]_{s})_{j\in\bb{\mathbf{d}}+1},
(\mathfrak{q}^{(2)}_{j})_{j\in\bb{\mathbf{d}}},
(U^{(1)}_{j})_{j\in\bb{\mathbf{d}}}
\right),
$$
such that 
$
\mathrm{sc}^{([1],2)}_{s}(\mathfrak{Q}^{(2)})=
\mathrm{tg}^{([1],2)}_{s}(\mathfrak{P}^{(2)}).
$

Then the \emph{$1$-composite}\index{composition!$1$-composition, $\circ^{1}$} of $\mathfrak{P}^{(2)}$ and $\mathfrak{Q}^{(2)}$, denoted by $\mathfrak{Q}^{(2)}\circ^{1}_{s}\mathfrak{P}^{(2)}$, is the ordered triple
$$
\mathfrak{Q}^{(2)}\circ^{1}_{s}\mathfrak{P}^{(2)}=
\left(
([R_{k}]_{s})_{k\in\bb{\mathbf{e}}+1},
(\mathfrak{r}^{(2)}_{k})_{k\in\bb{\mathbf{e}}},
(V^{(1)}_{k})_{k\in\bb{\mathbf{e}}}
\right),
$$
where $\mathbf{e}=\mathbf{c}\curlywedge\mathbf{d}$, and
\allowdisplaybreaks
\begin{alignat*}{2}
[R_{k}]_{s} &=
\begin{cases}
[P_{k}]_{s}, \\
[Q_{k-\bb{\mathbf{c}}}]_{s},    
\end{cases} 
&\qquad
&\begin{array}{l}
\text{if $k\in \bb{\mathbf{c}}+1$;} \\
\text{if $k\in [\bb{\mathbf{c}}+1,\bb{\mathbf{e}+1}]$,}
\end{array}
\\
\mathfrak{r}^{(2)}_{k} &=
\begin{cases}
\mathfrak{p}^{(2)}_{k},\\
\mathfrak{q}^{(2)}_{k-\bb{\mathbf{c}}},
\end{cases}
&\qquad
&\begin{array}{l}
\text{if $k\in \bb{\mathbf{c}}$;} \\
\text{if $k\in [\bb{\mathbf{c}},\bb{\mathbf{e}}]$,}
\end{array}
\\
V^{(1)}_{k} &=
\begin{cases}
T^{(1)}_{k},\\
U^{(1)}_{k-\bb{\mathbf{c}}},
\end{cases} 
&\qquad
&\begin{array}{l}
\text{if $k\in \bb{\mathbf{c}}$;} \\
\text{if $k\in [\bb{\mathbf{c}},\bb{\mathbf{e}}]$.}
\end{array}
\end{alignat*}
\end{restatable}

\begin{restatable}{proposition}{PDPthComp}
\label{PDPthComp}
Let $s$ be a sort in $S$, $\mathbf{c}, \mathbf{d}$ words in $S^{\star}$, $\mathfrak{P}^{(2)}$ a second-order $\mathbf{c}$-path in $\boldsymbol{\mathcal{A}}^{(2)}$, and $\mathfrak{Q}^{(2)}$ a second-order $\mathbf{d}$-path in $\boldsymbol{\mathcal{A}}^{(2)}$. Then, when defined, $\mathfrak{Q}^{(2)}\circ^{1}_{s}\mathfrak{P}^{(2)}$, the $1$-composite of $\mathfrak{P}^{(2)}$ and $\mathfrak{Q}^{(2)}$, is a second-order $(\mathbf{c}\curlywedge\mathbf{d})$-path in $\boldsymbol{\mathcal{A}}^{(2)}$ of the form
$$
\xymatrix@C=50pt{
\mathfrak{Q}^{(2)}\circ^{1}_{s}\mathfrak{P}^{(2)}
\colon
\mathrm{sc}^{([1],2)}_{s}(\mathfrak{P}^{(2)})
\ar@{=>}[r]^-{}
&
\mathrm{tg}^{([1],2)}_{s}(\mathfrak{Q}^{(2)}).
}
$$
Moreover, the above partial operation of $1$-composition, when defined, is associative, and, for every sort $s\in S$ and every path term class $[P]_{s}$ in $[\mathrm{PT}_{\boldsymbol{\mathcal{A}}}]_{s}$, the $(2,[1])$-identity second-order path on $[P]_{s}$ is, when defined, a neutral element for the operation of $1$-composition.
\end{restatable}

We next define the notion of subpath of a second-order path. But before doing that, since it will be used immediately below, we recall, from Definition~\ref{DSubw}, that, for a word $\mathbf{c}\in S^{\star}$ and indices $k$, $l\in\bb{\mathbf{c}}$ such that $k\leq l$, $\mathbf{c}^{k,l}$ is the subword  in $\mathbf{c}$ beginning at position $k$ and ending at position $l$.

\begin{restatable}{definition}{DDPthSub}
\label{DDPthSub}
Let $s$ be a sort in $S$, $\mathbf{c}$ a word in $S^{\star}$, $k$ and $l\in\bb{\mathbf{c}}$ such that $k\leq l$, and $\mathfrak{P}^{(2)}$ a second-order $\mathbf{c}$-path in $\mathrm{Pth}_{\boldsymbol{\mathcal{A}}^{(2)},s}$ of the form
$$\mathfrak{P}^{(2)}=\left(
([P_{i}]_{s})_{i\in\bb{\mathbf{c}}+1},
(\mathfrak{p}^{(2)}_{i})_{i\in\bb{\mathbf{c}}},
(T^{(1)}_{i})_{i\in\bb{\mathbf{c}}}
\right).
$$

Then we will denote by $\mathfrak{P}^{(2),k,l}$ the ordered triple
$$
\mathfrak{P}^{(2),k,l}
=
\left(
([P_{i}^{k,l}]_{s})_{i\in\bb{\mathbf{c}^{k,l}}+1},
(\mathfrak{p}^{(2),k,l}_{i})_{i\in\bb{\mathbf{c}^{k,l}}},
(T^{(1),k,l}_{i})_{i\in\bb{\mathbf{c}^{k,l}}}
\right)
$$
where, 
\begin{enumerate}
\item for every $i\in \bb{\mathbf{c}^{k,l}}+1$,
\begin{multicols}{2}
\noindent
\begin{enumerate}
\item[(i)] $[P_{i}^{k,l}]_{s}=[P_{k+i}]_{s}$;
\end{enumerate}
\end{multicols}
\item for every  $i\in \bb{\mathbf{c}^{k,l}}$,
\begin{multicols}{3}
\begin{enumerate}
\item[(i)]  $\mathfrak{p}^{(2), k,l}_{i}=\mathfrak{p}^{(2)}_{k+i}$;
\item[(ii)] $T^{(1),k,l}_{i}=T^{(1)}_{k+i}$.
\end{enumerate}
\end{multicols}
\end{enumerate}
We will call $\mathfrak{P}^{(2),k,l}$ the \emph{subpath of $\mathfrak{P}^{(2)}$ beginning at position $k$ and ending at position $l+1$}\index{path!second-order!subpath}. In particular, subpaths of the form  $\mathfrak{P}^{(2),0,k}$ will be called \emph{initial subpaths of $\mathfrak{P}^{(2)}$}, and subpaths of the form  $\mathfrak{P}^{(2), l,\bb{\mathbf{c}}-1}$ will be called \emph{final subpaths of $\mathfrak{P}^{(2)}$}.
\end{restatable}

In the following proposition we prove that the subpaths of a second-order path are actually second-order paths.

\begin{restatable}{proposition}{PDPthSub}
\label{PDPthSub}
Let $s$ be a sort in $S$, $\mathbf{c}$ be a word in $S^{\star}$, $k$ and $l$ indices in $\bb{\mathbf{c}}$ such that $k\leq l$, and $\mathfrak{P}^{(2)}$ a second-order $\mathbf{c}$-path in $\mathrm{Pth}_{\boldsymbol{\mathcal{A}}^{(2)},s}$ of the form
$$\mathfrak{P}^{(2)}=\left(
([P_{i}]_{s})_{i\in\bb{\mathbf{c}}+1},
(\mathfrak{p}^{(2)}_{i})_{i\in\bb{\mathbf{c}}},
(T^{(1)}_{i})_{i\in\bb{\mathbf{c}}}
\right).
$$

Then $\mathfrak{P}^{(2), k,l}$ 
is a second-order $\mathbf{c}^{k,l}$-path in $\mathrm{Pth}_{\boldsymbol{\mathcal{A}}^{(2)},s}$ of the form
$$
\xymatrix@C=50pt{
\mathfrak{P}^{(2), k,l}
\colon
[P_{k}]_{s}
\ar@{=>}[r]^-{}
&
[P_{l+1}]_{s}.
}
$$
\end{restatable}
\begin{proof}
The proof is in all respects similar to that of Proposition~\ref{PPthSub}.
\end{proof}

We next describe how the notion of subpath of a second-order path relates to that of $1$-composition of second-order paths.

\begin{restatable}{proposition}{PDPthRecons}
\label{PDPthRecons}
Let $\mathbf{c}$ be a word in $S^{\star}$, $k$ and $l$ indices in $\bb{\mathbf{c}}$ such that $k\leq l$, and $\mathfrak{P}^{(2)}$
a second-order $\mathbf{c}$-path in $\boldsymbol{\mathcal{A}}^{(2)}$ of sort $s\in S$. Then
$$\mathfrak{P}^{(2)} = \mathfrak{P}^{(2),l+1,\bb{\mathbf{c}}-1}
\circ^{1}_{s} \mathfrak{P}^{(2), k,l}
\circ^{1}_{s}\mathfrak{P}^{(2),0,k-1}.$$
\end{restatable}

\section{
\texorpdfstring
{Basic results on second-order paths}
{Basic results}
}
In this section we explain how a second-order path carries out the transformation of a path term class and state the structural results concerning this process. In addition, we introduce the notion of second-order echelon, a key concept in the development of our theory.

\begin{restatable}{definition}{DDEch}
\label{DDEch} 
We denote by $\mathrm{ech}^{(2,\mathcal{A}^{(2)})}$ the $S$-sorted mapping
$$
\mathrm{ech}^{(2,\mathcal{A}^{(2)})}
\colon
\mathcal{A}^{(2)}
\mor
\mathrm{Pth}_{\mathcal{A}^{(2)}}
$$
which, for every sort $s\in S$, 
$$\textstyle
  \mathrm{ech}^{(2,\mathcal{A}^{(2)})}_{s}\nfunction
  {\mathcal{A}^{(2)}_{s}}
  {\mathrm{Pth}_{\boldsymbol{\mathcal{A}}^{(2)},s}}
  {\mathfrak{p}^{(2)}=([M]_{s},[N]_{s})}
  {\left(
\left([M]_{s},[N]_{s}\right),\mathfrak{p}^{(2)},
\mathrm{id}^{\mathbf{T}_{\Sigma^{\boldsymbol{\mathcal{A}}}}(X)_{s}}
\right)}
$$
%
%
%
%
We also let
$$
\xymatrix@C=90pt{
\mathrm{ech}^{(2,\mathcal{A}^{(2)})}_{s}(\mathfrak{p}^{(2)})\colon
[M]_{s}
\ar@{=>}[r]^-{\text{\Small{($\mathfrak{p}^{(2)}$, $
\mathrm{id}^{\mathbf{T}_{\Sigma^{\boldsymbol{\mathcal{A}}}}(X)_{s}}
$)}}}
&
[N]_{s}
}
$$
stand for $\mathrm{ech}^{(2,\mathcal{A}^{(2)})}_{s}(\mathfrak{p}^{(2)})$. This sequence, actually, satisfies the defining conditions of the notion of second-order path since (i) $\mathrm{id}^{\mathrm{T}_{\Sigma^{\boldsymbol{\mathcal{A}}}}(X)_{s}}(M)\in [M]_{s}$ and (ii) $\mathrm{id}^{\mathrm{T}_{\Sigma^{\boldsymbol{\mathcal{A}}}}(X)_{s}}(N)\in [N]_{s}$. We call $\mathrm{ech}^{(2,\mathcal{A}^{(2)})}_{s}(\mathfrak{p}^{(2)})$ the \emph{second-order echelon associated to}
$\mathfrak{p}^{(2)}$. 
\index{echelon!second-order!$\mathrm{ech}^{(2,\mathcal{A}^{(2)})}$}
Moreover, for a sort $s$ in $S$, we will say that a second-order path $\mathfrak{P}^{(2)}\in \mathrm{Pth}_{\boldsymbol{\mathcal{A}},s}$ is a \emph{second-order echelon} if there exists a second-order rewrite rule $\mathfrak{p}^{(2)}\in\mathcal{A}^{(2)}_{s}$ such that $\mathfrak{P}^{(2)} = \mathrm{ech}^{(2,\mathcal{A}^{(2)})}_{s}(\mathfrak{p}^{(2)})$, i.e., if $\mathfrak{P}^{(2)}\in \mathrm{Im}(\mathrm{ech}^{(2,\mathcal{A}^{(2)})})_{s}$.  We denote by $\mathrm{Ech}^{(2,\mathcal{A}^{(2)})}[\mathcal{A}^{(2)}]$ the $S$-sorted set $\mathrm{Im}(\mathrm{ech}^{(2,\mathcal{A}^{(2)})})$. Finally, we will say that a second-order path $\mathfrak{P}^{(2)}$ is \emph{echelonless}\index{path!second-order!echelonless} if $\bb{\mathfrak{P}^{(2)}}\geq 1$ and none of its one-step subpaths is a second-order echelon. 
\end{restatable}

From the above it follows, as stated in the following corollary, that the first-order translations of an echelonless second-order path must be non-identity translations.

\begin{restatable}{corollary}{CDEch}
\label{CDEch}
Let $s$ be a sort in $S$, $\mathbf{c}\in S^{\star}$ and $\mathfrak{P}^{(2)}$ an echelonless second-order $\mathbf{c}$-path in $\mathrm{Pth}_{\boldsymbol{\mathcal{A}}^{(2)},s}$ of the form
$$
\mathfrak{P}^{(2)}
=
\left(
(
[P_{i}]_{s}
)_{i\in\bb{\mathbf{c}}+1}, 
(
\mathfrak{p}^{(2)}_{i}
)_{i\in\bb{\mathbf{c}}},
(
T^{(1)}_{i}
)_{i\in\bb{\mathbf{c}}}
\right).
$$

Then $\mathbf{c}\neq \lambda$ and, for every $i\in \bb{\mathbf{c}}$, $T^{(1)}_{i}\neq \mathrm{id}^{\mathrm{T}_{\Sigma^{\boldsymbol{\mathcal{A}}}}(X)_{s}}$.
\end{restatable}

In the following definition we introduce the notion of a head-constant echelonless second-order path.

\begin{restatable}{definition}{DDHeadCt}
\label{DDHeadCt}
Let $s$ be a sort in $S$, $\mathbf{c}$ a word in $S^{\star}-\{\lambda\}$ and  $\mathfrak{P}^{(2)}$ an echelonless second-order $\mathbf{c}$-path in $\mathrm{Pth}_{\mathcal{A}^{(2)},s}$ of the form 
$$
\mathfrak{P}^{(2)}
=
\left(
(
[P_{i}]_{s}
)_{i\in\bb{\mathbf{c}}+1}, 
(
\mathfrak{p}^{(2)}_{i}
)_{i\in\bb{\mathbf{c}}},
(
T^{(1)}_{i}
)_{i\in\bb{\mathbf{c}}}
\right).
$$
We will say that $\mathfrak{P}^{(2)}$ is a \emph{head-constant} echelonless second-order path\index{path!second-order!head-constant}\index{head-constant} if $(T^{(1)}_{i})_{i\in\bb{\mathbf{c}}}$, the family of the first-order translations occurring in it, have the same type, that is, they are associated to the same operation symbol.  
\end{restatable}

Unlike in the first part of this work, an echelonless second-order path does not necessarily goes across head-constant families of first-order translations. This is shown in the following example.

\begin{example}\label{EHeadCt} For the set of sorts $S=1$, the  signature $\Sigma=\{\sigma\}$, where $\sigma$ is a binary operation, and for $X=\{x,y,z\}$, we consider the free $\Sigma$-algebra on $X$, that is $\mathbf{T}_{\Sigma}(X)$. In the example, we will omit any reference to the unique sort in $S$. Consider the set of rewrite rules $\mathcal{A}=\{\mathfrak{p},\mathfrak{q},\mathfrak{r}\}$, where
\begin{align*}
\mathfrak{p}&=(x,y),&
\mathfrak{q}&=(y,z),&
\mathfrak{r}&=(x,z).
\end{align*}
and the set of second-order rewrite rules $\mathcal{A}^{(2)}=\{\mathfrak{p}^{(2)}\}$, where 
$$\mathfrak{p}^{(2)}=([\mathfrak{q}\circ^{0}\mathfrak{p}],[\mathfrak{r}]).$$

Let us note that $\mathfrak{p}^{(2)}$ is a valid second-order rewrite rule because
\allowdisplaybreaks
\begin{align*}
\mathrm{sc}^{(0,1)}(\mathrm{ip}^{(1,X)@}(\mathfrak{q}\circ^{0}\mathfrak{p}))
&=
\mathrm{sc}^{(0,1)}(\mathrm{ip}^{(1,X)@}(\mathfrak{r}))
=
x;
\\
\mathrm{tg}^{(0,1)}(\mathrm{ip}^{(1,X)@}(\mathfrak{q}\circ^{0}\mathfrak{p}))
&=
\mathrm{tg}^{(0,1)}(\mathrm{ip}^{(1,X)@}(\mathfrak{r}))
=z.
\end{align*}

Consider the second-order path
$$
\xymatrix@C=150pt@R=15pt{
\mathfrak{P}^{(2)}: [
\sigma(\mathfrak{q}\circ^{0}\mathfrak{p}, \mathfrak{q}\circ^{0}\mathfrak{p})
]
\ar@2{->}[r]^-{(\mathfrak{p}^{(2)}, \sigma(\underline{\quad}, \mathfrak{q}\circ^{0}\mathfrak{p}))}
&
[\sigma(\mathfrak{r}, \mathfrak{q}\circ^{0}\mathfrak{p})]
\\
{\qquad\qquad\qquad\qquad\qquad\quad}
\ar@2{->}[r]^-{(\mathfrak{p}^{(2)}, 
\sigma(\mathfrak{r},z) \circ^{0} \sigma(x, \underline{\quad}))}
&
\quad
[\sigma(\mathfrak{r},\mathfrak{r})].
\quad
}
$$

One can easily verify that  $\mathfrak{P}^{(2)}$ is a second-order path from $[\sigma(\mathfrak{q}\circ^{0}\mathfrak{p}, \mathfrak{q}\circ^{0}\mathfrak{p})]$ to $[\sigma(\mathfrak{r},\mathfrak{r})]$ of length $2$.  Let us note that none of the one-step subpaths of $\mathfrak{P}^{(2)}$ is a  second-order echelon. Thus $\mathfrak{P}^{(2)}$ is an echelonless second-order path.

Therefore, in contrast to what happens in the first part of this work, the translations $\sigma(\mathfrak{r}, \underline{\quad})$ and $\sigma(\mathfrak{r},z) \circ^{0} \sigma(x, \underline{\quad})$ are non-identity translations but they are not associated to the same operation symbol.
\end{example}

We next introduce the notion of coherent head-constant echelonless second-order path.

\begin{restatable}{definition}{DDCoh}
\label{DDCoh}\index{path!second-order!coherent}\index{coherent} Let $s$ be a sort in $S$, $\mathbf{c}$ a word in $S^{\star}-\{\lambda\}$ and  $\mathfrak{P}^{(2)}$ a head-constant echelonless second-order $\mathbf{c}$-path in $\mathrm{Pth}_{\mathcal{A}^{(2)},s}$ of the form 
$$
\mathfrak{P}^{(2)}
=
\left(
(
[P_{i}]_{s}
)_{i\in\bb{\mathbf{c}}+1}, 
(
\mathfrak{p}^{(2)}_{i}
)_{i\in\bb{\mathbf{c}}},
(
T^{(1)}_{i}
)_{i\in\bb{\mathbf{c}}}
\right)
$$
where, for a unique word $\mathbf{s}\in S^{\star}-\{\lambda\}$, a unique operation symbol $\tau\in \Sigma^{\boldsymbol{\mathcal{A}}}_{\mathbf{s},s}$, the family  $(
T^{(1)}_{i}
)_{i\in\bb{\mathbf{c}}}$ is a family of first-order translations of type $\tau$.
That is, for every $i\in \bb{\mathbf{c}}$, there exists a unique index $k_{i}\in \bb{\mathbf{s}}$, a unique family of paths $(\mathfrak{P}_{i,j})_{j\in k}\in \prod_{j\in k}\mathrm{Pth}_{\boldsymbol{\mathcal{A}},s_{j}}$, a unique family of paths $(\mathfrak{P}_{i,l})_{l\in \bb{\mathbf{s}}-(k+1)}\in \prod_{l\in \bb{\mathbf{s}}-(k+1)}\mathrm{Pth}_{\boldsymbol{\mathcal{A}},s_{l}}$ and a unique first-order translation $T^{(1)'}_{i}\in \mathrm{Tl}_{c_{i}}(\mathbf{T}_{\Sigma^{\boldsymbol{\mathcal{A}}}}(X))_{s_{k}}$  and
\allowdisplaybreaks
\begin{multline*}
T^{(1)}_{i}=\tau^{\mathbf{T}_{\Sigma^{\boldsymbol{\mathcal{A}}}}(X)}\left(
\mathrm{CH}^{(1)}_{s_{0}}\left(\mathfrak{P}_{i,0}\right),
\cdots,
\mathrm{CH}^{(1)}_{s_{k_{i}-1}}\left(\mathfrak{P}_{i,k_{i}-1}\right),
\right.
\\
\left.
T^{(1)'}_{i},
\mathrm{CH}^{(1)}_{s_{k_{i}+1}}\left(\mathfrak{P}_{i,k_{i}+1}\right),
\cdots,
\mathrm{CH}^{(1)}_{s_{\bb{\mathbf{s}}-1}}\left(\mathfrak{P}_{i,\bb{\mathbf{s}}-1}\right)
\right).
\end{multline*}
Let us assume that, for every $i\in\bb{\mathbf{c}}$, the second-order rewrite rule $\mathfrak{p}^{(2)}_{i}$ is given by $([M_{i}]_{c_{i}}, [N_{i}]_{c_{i}})$. Since $\mathfrak{P}^{(2)}$ is a second-order path, we have that, for every $i\in\bb{\mathbf{c}}$, it is the case that 
\begin{itemize}
\item[(i)] $T^{(1)}_{i}(M_{i})\in [P_{i}]_{s}$;
\item[(ii)] $T^{(1)}_{i}(N_{i})\in [P_{i+1}]_{s}$.
\end{itemize}
In particular, for every $i\in\bb{\mathbf{c}}-1$, the following equality holds
\[
\left[
T^{(1)}_{i}\left(
N_{i}\right)
\right]_{s}
=
\left[
T^{(1)}_{i+1}\left(
M_{i+1}\right)
\right]_{s}.
\]
That is, the following equality holds
\allowdisplaybreaks
\begin{multline*}
\left[
\tau^{\mathbf{T}_{\Sigma^{\boldsymbol{\mathcal{A}}}}(X)}\left(
\mathrm{CH}^{(1)}_{s_{0}}\left(\mathfrak{P}_{i,0}\right),
\cdots,
\mathrm{CH}^{(1)}_{s_{k_{i}-1}}\left(\mathfrak{P}_{i,k_{i}-1}\right),
\right.\right.
\\
\left.\left.
T^{(1)'}_{i}\left(
N_{i}
\right)
,
\mathrm{CH}^{(1)}_{s_{k_{i}+1}}\left(\mathfrak{P}_{i,k_{i}+1}\right),
\cdots,
\mathrm{CH}^{(1)}_{s_{\bb{\mathbf{s}}-1}}\left(\mathfrak{P}_{i,\bb{\mathbf{s}}-1}\right)
\right)
\right]_{s}
\\=
\left[
\tau^{\mathbf{T}_{\Sigma^{\boldsymbol{\mathcal{A}}}}(X)}\left(
\mathrm{CH}^{(1)}_{s_{0}}\left(\mathfrak{P}_{i+1,0}\right),
\cdots,
\mathrm{CH}^{(1)}_{s_{k_{i+1}-1}}\left(\mathfrak{P}_{i+1,k_{i+1}-1}\right),
\right.\right.
\\
\left.\left.
T^{(1)'}_{i+1}\left(
M_{i+1}
\right)
,
\mathrm{CH}^{(1)}_{s_{k_{i+1}+1}}\left(\mathfrak{P}_{i,k_{i+1}+1}\right),
\cdots,
\mathrm{CH}^{(1)}_{s_{\bb{\mathbf{s}}-1}}\left(\mathfrak{P}_{i+1,\bb{\mathbf{s}}-1}\right)
\right)
\right]_{s}.
\end{multline*}
We will say that $\mathfrak{P}^{(2)}$ is \emph{coherent} if, for every $i\in\bb{\mathbf{c}}-1$, from the above equality, we can derive the following $\bb{\mathbf{s}}$ equalities
\allowdisplaybreaks
\begin{align*}
\left[
\mathrm{CH}^{(1)}_{s_{0}}\left(\mathfrak{P}_{i,0}\right)
\right]_{s_{0}}
&=
\left[
\mathrm{CH}^{(1)}_{s_{0}}\left(\mathfrak{P}_{i+1,0}\right)
\right]_{s_{0}}
\\
&\,\,\,\,\vdots
\\
\left[
T^{(1)'}_{i}\left(
N_{i}
\right)
\right]_{s_{k_{i}}}
&=
\left[
\mathrm{CH}^{(1)}_{s_{k_{i}}}\left(\mathfrak{P}_{i+1,k_{i}}\right)
\right]_{s_{k_{i}}}
\\
&\,\,\,\,\vdots
\\
\left[
\mathrm{CH}^{(1)}_{s_{k_{i+1}}}\left(\mathfrak{P}_{i,k_{i+1}}\right)
\right]_{s_{k_{i+1}}}
&=
\left[
T^{(1)'}_{i+1}\left(
M_{i+1}
\right)
\right]_{s_{k_{i+1}}}
\\
&\,\,\,\,\vdots
\\
\left[
\mathrm{CH}^{(1)}_{s_{\bb{\mathbf{s}}-1}}\left(\mathfrak{P}_{i,\bb{\mathbf{s}}-1}\right)
\right]_{s_{\bb{\mathbf{s}}-1}}
&=
\left[
\mathrm{CH}^{(1)}_{s_{\bb{\mathbf{s}}-1}}\left(\mathfrak{P}_{i+1,\bb{\mathbf{s}}-1}\right)
\right]_{s_{\bb{\mathbf{s}}-1}}.
\end{align*}

We warn the reader that, in the above $\bb{\mathbf{s}}$ equalities, we have implicitly assumed that $k_{i}<k_{i+1}$, however the cases $k_{i}=k_{i+1}$ and $k_{i+1}<k_{i}$ should also be taken into account. We have not written them down so as not to lengthen the presentation of this notion.
\end{restatable}

We next present examples of second-order paths illustrating the aforementioned properties.

\begin{example}\label{EDLoop}[Continuation of Example~\ref{EHeadCt}] 
Consider the same second-order rewrite system as in Example~\ref{EHeadCt}.
%
Consider the second-order path
$$
\xymatrix@C=150pt@R=15pt{
\mathfrak{Q}^{(2)}: [
\sigma(\mathfrak{q}\circ^{0}\mathfrak{p}, \mathfrak{q}\circ^{0}\mathfrak{p})
]
\ar@2{->}[r]^-{(\mathfrak{p}^{(2)}, 
\sigma(\underline{\quad},z)
\circ^{0}
\sigma(x,\mathfrak{q}\circ^{0}\mathfrak{p}))}
&
[\sigma(\mathfrak{r}, \mathfrak{q}\circ^{0}\mathfrak{p})]
\\
{\qquad\qquad\qquad\qquad\qquad\quad}
\ar@2{->}[r]^-{(\mathfrak{p}^{(2)}, 
\sigma(z,z) \circ^{0} \sigma(\mathfrak{r}, \underline{\quad}))}
&
\quad
[\sigma(\mathfrak{r},\mathfrak{r})].
\quad
}
$$

One can easily verify that  $\mathfrak{Q}^{(2)}$ is a head-constant echelonless second-order path from $[\sigma(\mathfrak{q}\circ^{0}\mathfrak{p}, \mathfrak{q}\circ^{0}\mathfrak{p})]$ to $[\sigma(\mathfrak{r},\mathfrak{r})]$ of length $2$ associated to the $0$-composition operation symbol $\circ^{0}$. However, this second-order path is not coherent. This is so, because from the equality
\[
\left(
\left[
\sigma(\mathfrak{r},z)
\right]
\right)
\circ^{0[\mathbf{PT}_{\boldsymbol{\mathcal{A}}}]}
\left(
\left[
\sigma(x,\mathfrak{q}\circ^{0}\mathfrak{p}
\right]
\right)
=
\left(
\left[
\sigma(z,z)
\right]
\right)
\circ^{0[\mathbf{PT}_{\boldsymbol{\mathcal{A}}}]}
\left(
\left[
\sigma(\mathfrak{r},\mathfrak{q}\circ^{0}\mathfrak{p}
\right]
\right),
\]
we cannot infer that 
$
[\sigma(\mathfrak{r},z)]
=
[\sigma(z,z)].
$ In this regard, let us note that the class on the left represents paths of length $1$, while the class on the right represents paths of length $0$.

Now consider the second-order path
$$
\xymatrix@C=150pt@R=15pt{
\mathfrak{R}^{(2)}: [
\sigma(\mathfrak{q}\circ^{0}\mathfrak{p}, \mathfrak{q}\circ^{0}\mathfrak{p})
]
\ar@2{->}[r]^-{(\mathfrak{p}^{(2)}, 
\sigma(\underline{\quad}, z)
\circ^{0}
\sigma(x,\mathfrak{q}\circ^{0}\mathfrak{p}))}
&
[\sigma(\mathfrak{r}, \mathfrak{q}\circ^{0}\mathfrak{p})]
\\
{\qquad\qquad\qquad\qquad\qquad\quad}
\ar@2{->}[r]^-{(\mathfrak{p}^{(2)}, 
\sigma(\mathfrak{r},z) \circ^{0} \sigma(x, \underline{\quad}))}
&
\quad
[\sigma(\mathfrak{r},\mathfrak{r})].
\quad
}
$$
This is also a  head-constant echelonless second-order path associated to the $0$-composition operation symbol $\circ^{0}$. However, this second-order path is coherent. 
\end{example}

From the definitions of the notions of coherent and head-constant echelonless second-order paths, introduced in Definitions~\ref{DDHeadCt} and~\ref{DDCoh}, respectively, we can immediately infer that the echelonless one-step second-order paths are always head-constant and coherent.

\begin{restatable}{corollary}{CDCohOneStep}
\label{CDCohOneStep} 
Let $s$ be a sort in $S$ and  $\mathfrak{P}^{(2)}$ a one-step echelonless second-order path in $\mathrm{Pth}_{\mathcal{A}^{(2)},s}$. Then $\mathfrak{P}^{(2)}$ is a coherent head-constant echelonless second-order path.
\end{restatable}
\begin{proof}
Since $\bb{\mathfrak{P}^{(2)}}=1$, we have that the definitions of head-constant and coherent are trivially fulfilled.
\end{proof}

In the following corollary we state that the head-constant echelonless second-order paths that are associated to an operation symbol coming from $\Sigma$, the original signature, are coherent.

\begin{restatable}{corollary}{CDCohSigma}
\label{CDCohSigma} Let $s$ be a sort in $S$, $\mathbf{c}$ a word in $S^{\star}-\{\lambda\}$ and  $\mathfrak{P}^{(2)}$ a head-constant echelonless second-order $\mathbf{c}$-path in $\mathrm{Pth}_{\mathcal{A}^{(2)},s}$ of the form 
$$
\mathfrak{P}^{(2)}
=
\left(
(
[P_{i}]_{s}
)_{i\in\bb{\mathbf{c}}+1}, 
(
\mathfrak{p}^{(2)}_{i}
)_{i\in\bb{\mathbf{c}}},
(
T^{(1)}_{i}
)_{i\in\bb{\mathbf{c}}}
\right)
$$
where, for a unique $\mathbf{s}\in S^{\star}$, $(
T^{(1)}_{i}
)_{i\in\bb{\mathbf{c}}}$ is a family of first-order translations of type associated  to a unique operation symbol $\sigma\in \Sigma_{\mathbf{s},s}$.
Then $\mathfrak{P}^{(2)}$ is coherent.
\end{restatable}
\begin{proof}
It follows from Lemma~\ref{LPTQSigma}.
\end{proof}

The following lemma is fundamental for the development of the present work. In it we show that for every 
coherent head-constant echelonless second-order path $\mathfrak{P}^{(2)}$ associated to an operation symbol $\tau\in \Sigma^{\boldsymbol{\mathcal{A}}}_{\mathbf{s},s}$ has canonically associated a family $(\mathfrak{P}^{(2)}_{j})_{j\in\bb{\mathbf{s}}}$ of second-order paths such that, for every $j\in\bb{\mathbf{s}}$, $\mathfrak{P}^{(2)}_{j}$ is a second-order path of sort $s_{j}$ and they are in charge of transforming in parallel each of the components of $\mathfrak{P}^{(2)}$.

\begin{restatable}{lemma}{LDPthExtract}
\label{LDPthExtract}\index{path!second-order!extraction procedure}
Let $s$ be a sort in $S$, $\mathbf{c}$ a word in $S^{\star}-\{\lambda\}$, and $\mathfrak{P}^{(2)}$
a coherent head-constant echelonless second-order $\mathbf{c}$-path in $\mathrm{Pth}_{\mathcal{A}^{(2)},s}$ of the form 
$$
\mathfrak{P}^{(2)}=\left(
([P_{i}]_{s})_{i\in\bb{\mathbf{c}}+1},
(\mathfrak{p}^{(2)}_{i})_{i\in\bb{\mathbf{c}}},
(T^{(1)}_{i})_{i\in\bb{\mathbf{c}}}
\right)
$$
where, for a unique $\mathbf{s}\in S^{\star}$, $(T^{(1)}_{i})_{i\in\bb{\mathbf{c}}}$ is a family of first-order translations of type $\tau$, for a unique operation symbol $\tau\in \Sigma^{\boldsymbol{\mathcal{A}}}_{\mathbf{s},s}$. Then there exists a unique pair
$
\textstyle
((\mathbf{c}_{j})_{j\in\bb{\mathbf{s}}},(\mathfrak{P}^{(2)}_{j})_{j\in\bb{\mathbf{s}}})\in (S^{\star})^{\bb{\mathbf{s}}}\times \mathrm{Pth}_{\boldsymbol{\mathcal{A}}^{(2)},\mathbf{s}}
$
such that, for every $j\in\bb{\mathbf{s}}$, $\mathfrak{P}^{(2)}_{j}$ is a second-order $\mathbf{c}_{j}$-path in $\mathrm{Pth}_{\mathcal{A}^{(2)},s_{j}}$ of the form
$$
\mathfrak{P}^{(2)}_{j}=
\left(
([P_{j,k}]_{s_{j}})_{k\in\bb{\mathbf{c}_{j}}+1},
(\mathfrak{p}^{(2)}_{j,k})_{k\in\bb{\mathbf{c}_{j}}},
(T^{(1)}_{j,k})_{k\in\bb{\mathbf{c}_{j}}}
\right)$$
and there exists a unique mapping $(i^{j,k})_{(j,k)\in \coprod_{j\in \bb{\mathbf{s}}}\bb{\mathbf{c}_{j}}}\colon \coprod_{j\in \bb{\mathbf{s}}}\bb{\mathbf{c}_{j}}\mor \bb{\mathbf{c}}$ such that, for every $(j,k)$ in $\coprod_{j\in \bb{\mathbf{s}}}\bb{\mathbf{c}_{j}}$, the following equalities 
\begin{enumerate}
\item[(i)] $c_{j,k}=c_{i^{j,k}}$ and
\item[(ii)] $\mathfrak{p}^{(2)}_{j,k}=\mathfrak{p}^{(2)}_{i^{j,k}}$
\end{enumerate}
are fulfilled.
\end{restatable}
\begin{proof}
Let us assume that, for every $i\in\bb{\mathbf{c}}$, the second-order rewrite rule $\mathfrak{p}^{(2)}_{i}$ is given by $([M_{i}]_{c_{i}}, [N_{i}]_{c_{i}})$. Since $\mathfrak{P}^{(2)}$ is a second-order path, we have that, for every $i\in\bb{\mathbf{c}}$, it is the case that 
\begin{itemize}
\item[(i)] $T^{(1)}_{i}\left(
M_{i}
\right)\in [P_{i}]_{s}$ and
\item[(ii)] $T^{(1)}_{i}\left(
N_{i}
\right)\in [P_{i+1}]_{s}$.
\end{itemize}

Let us note that, since $\mathfrak{P}^{(2)}$ is echelonless, for every $i\in\bb{\mathbf{c}}$, we have that $T^{(1)}_{i}\neq \mathrm{id}^{\mathbf{T}_{\Sigma^{\boldsymbol{\mathcal{A}}}}(X)_{s}}$. Since the family of first-order translations $(
T^{(1)}_{i}
)_{i\in\bb{\mathbf{c}}}$ has type $\tau\in \Sigma^{\boldsymbol{\mathcal{A}}}_{\mathbf{s},s}$, for every $i\in\bb{\mathbf{c}}$, there exists a unique index $j_{i}\in \bb{\mathbf{s}}$, a unique family of paths $(\mathfrak{P}_{i,k})_{k\in j_{i}}\in \prod_{k\in j_{i}}\mathrm{Pth}_{\boldsymbol{\mathcal{A}},s_{k}}$, a unique family of paths $(\mathfrak{P}_{i,l})_{l\in \bb{\mathbf{s}}-(j_{i}+1)}\in \prod_{l\in \bb{\mathbf{s}}-(j_{i}+1)}\mathrm{Pth}_{\boldsymbol{\mathcal{A}},s_{l}}$ and a unique first-order translation $T^{(1)'}_{i}\in \mathrm{Tl}_{c_{i}}(\mathbf{T}_{\Sigma^{\boldsymbol{\mathcal{A}}}}(X))_{s_{j_{i}}}$  such that 
\begin{multline*}
T^{(1)}_{i}=\tau^{\mathbf{T}_{\Sigma^{\boldsymbol{\mathcal{A}}}}(X)}\left(
\mathrm{CH}^{(1)}_{s_{0}}\left(\mathfrak{P}_{i,0}\right),
\cdots,
\mathrm{CH}^{(1)}_{s_{j_{i}-1}}\left(\mathfrak{P}_{i,j_{i}-1}\right),
\right.
\\
\left.
T^{(1)'}_{i},
\mathrm{CH}^{(1)}_{s_{j_{i}+1}}\left(\mathfrak{P}_{i,j_{i}+1}\right),
\cdots,
\mathrm{CH}^{(1)}_{s_{\bb{\mathbf{s}}-1}}\left(\mathfrak{P}_{i,\bb{\mathbf{s}}-1}\right)
\right).
\end{multline*}
In this case we will say that $\mathfrak{P}^{(2),i,i}$ is a one-step second-order $c_{i}$-path of type $j_{i}$.

For every $j\in\bb{\mathbf{s}}$, let $\mathbf{c}_{j}$ be the word in $S^{\star}$ obtained by concatenating,  in order of appearance, the letters $c_{i}\in S$ for which $i\in\bb{\mathbf{c}}$ is such that $\mathfrak{P}^{(2),i,i}$ is a one-step second-order $c_{i}$-path of type $j$.  

If $\mathbf{c}_{j}=\lambda$, i.e., if there are no steps of type $j$, we define
\[\mathfrak{P}^{(2)}_{j}=\mathrm{ip}^{(2,[1])\sharp}_{s_{j}}\left(
\left[
\mathrm{CH}^{(1)}_{s_{j}}\left(
\mathfrak{P}_{0,j}
\right)
\right]_{s_{j}}
\right).
\]

Let us note that the coherence condition on $\mathfrak{P}^{(2)}$ implies that the election of the specific step does not matter, as it happens that, for every $i\in\bb{\mathbf{c}}$,
\[
\left[
\mathrm{CH}^{(1)}_{s_{j}}\left(
\mathfrak{P}_{0,j}
\right)
\right]_{s_{j}}
=
\left[
\mathrm{CH}^{(1)}_{s_{j}}\left(
\mathfrak{P}_{i,j}
\right)
\right]_{s_{j}}
.
\]

We will therefore focus our attention on the case in which $\mathbf{c}_{j}\neq \lambda$. 
By construction, for every $j\in\bb{\mathbf{s}}$ and every $k\in\bb{\mathbf{c}_{j}}$, there exists a unique $i^{j,k}\in\bb{\mathbf{c}}$ such that $\mathfrak{P}^{(2),i^{j,k},i^{j,k}}$ is the $k$-th one-step second-order path of type $j$, which starts at position $i^{j,k}$ and ends at position $i^{j,k}+1$. Hence, $c_{j,k}=c_{i^{j,k}}$, i.e., the $k$-th letter in the word $\mathbf{c}_{j}$ is equal to the letter $c_{i^{j,k}}$ in the word $\mathbf{c}$.

Now, for every $j\in\bb{\mathbf{s}}$, let $\mathfrak{P}^{(2)}_{j}$ be the ordered triple  
$$\mathfrak{P}^{(2)}_{j} = \left(([P_{j,k}]_{s_{j}})_{k\in\bb{\mathbf{c}_{j}}+1},  (\mathfrak{p}^{(2)}_{j, k})_{k\in\bb{\mathbf{c}_{j}}},
(T^{(1)}_{ j, k})_{k\in\bb{\mathbf{c}_{j}}}
\right)$$ 
defined as follows:
\allowdisplaybreaks
\begin{alignat*}{2}
P_{j,k} &=
\begin{cases}
T^{(1)'}_{i^{j,k}}\left(
M_{i^{j,k}}
\right)
, \\
T^{(1)'}_{i^{j,\bb{\mathbf{c}_{j}}-1}}\left(
N_{{i^{j,\bb{\mathbf{c}_{j}}-1}}}
\right),    
\end{cases} 
&\qquad
&\begin{array}{l}
\text{if $k\in \bb{\mathbf{c}_{j}}$;} \\
\text{if $k=\bb{\mathbf{c}_{j}}$;}
\end{array}
\\
\mathfrak{p}^{(2)}_{j,k} &=
\begin{cases}
\mathfrak{p}^{(2)}_{i^{j,k}},   
\end{cases} 
&\qquad
&\begin{array}{l}
\text{if $k\in \bb{\mathbf{c}_{j}}$;} 
\end{array}
\\
T^{(1)}_{j,k} &=
\begin{cases}
T^{(1)'}_{i^{j,k}},   
\end{cases} 
&\qquad
&\begin{array}{l}
\text{if $k\in \bb{\mathbf{c}_{j}}$.} 
\end{array}
\end{alignat*}
That is,  (1)  for every $k\in\bb{\mathbf{c}_{j}}$, $P_{j,k}$, i.e., the $k$-th path term representative  of the sequence $([P_{ j, k}]_{s_{j}})_{k\in\bb{\mathbf{c}_{j}}+1}$ is equal to the $j$-th subterm of the path term $T_{i^{j,k}}(
M_{i^{j,k}}
)$, i.e., $T^{(1)'}_{i^{j,k}}(
M_{i^{j,k}}
)$; (2) for every $k\in\bb{\mathbf{c}_{j}}$, $\mathfrak{p}^{(2)}_{j,k}$ is equal to $\mathfrak{p}^{(2)}_{i^{j,k}}$, i.e., the $k$-th second-order rewrite rule of the sequence $(\mathfrak{p}^{(2)}_{j,k})_{k\in\bb{\mathbf{c}_{j}}}$ is equal to the second-order rewrite rule associated to the $k$-th index of type $j$; and, finally, (3) for every $k\in\bb{\mathbf{c}_{j}}$,  $T^{(1)}_{ j, k}$ is equal to $T^{(1)'}_{i^{j,k}}$, i.e., the $k$-th first-order translation of the sequence $(T^{(1)}_{ j, k})_{k\in\bb{\mathbf{c}_{j}}}$ is equal to $T^{(1)'}_{i^{j,k}}$ the derived first-order translation associated to the $k$-th first-order translation of type $j$.

We claim that, for every $j\in\bb{\mathbf{s}}$, $\mathfrak{P}^{(2)}_{j}$ is a second-order $\mathbf{c}_{j}$-path of sort $s_{j}$ in $\boldsymbol{\mathcal{A}}^{(2)}$ from the path term class
$[
T^{(1)'}_{i^{j,0}}(
M_{i^{j,0}}
)]_{s_{j}}$ to the path term class 
$[
T^{(1)'}_{i^{j,\bb{\mathbf{c}_{j}}-1}}(
N_{{i^{j,\bb{\mathbf{c}_{j}}-1}}}
)]_{s_{j}}$. These are, respectively, the classes associated to the first and last path terms in the sequence $(P_{j,k})_{k\in\bb{\mathbf{c}_{j}}+1}$.

Let us note that, for every $k\in\bb{\mathbf{c}_{j}}$, we have that $\mathfrak{p}^{(2)}_{j, k}=\mathfrak{p}^{(2)}_{i^{j,k}}$, where $\mathfrak{p}^{(2)}_{i^{j,k}}=([M_{i^{j,k}}]_{c_{i^{j,k}}}, [N_{i^{j,k}}]_{c_{i^{j,k}}})$, and $T^{(1)}_{j,k}=T^{(1)'}_{i^{j,k}}$. We have that, for every $k\in\bb{\mathbf{c}_{j}}$, the following condition immediately follows from the definition of $P_{j,k}$
\begin{enumerate}
\item[(i)] $T^{(1)'}_{i^{j,k}}\left(
M_{i^{j,k}}
\right)
\in [P_{j,k}]_{s_{j}}.
$
\end{enumerate}

It is left to us to show that, for every $k\in\bb{\mathbf{c}_{j}}$, the following condition also holds
\begin{enumerate}
\item[(ii)] $T^{(1)'}_{i^{j,k}}\left(
N_{i^{j,k}}
\right)
\in [P_{j,k+1}]_{s_{j}}.
$
\end{enumerate}

Let us note that for $k=\bb{\mathbf{c}_{j}}-1$, the condition (ii) above follows, immediately, by the definition of $P_{j,\bb{\mathbf{c}_{j}}}$. Therefore, it suffices to prove condition (ii) above for every $k\in\bb{\mathbf{c}_{j}}-1$. That is, we need to prove that, for every $k\in\bb{\mathbf{c}_{j}}-1$, the following condition holds 
\[
T^{(1)'}_{i^{j,k}}\left(
N_{i^{j,k}}
\right)
\in 
\left[
T^{(1)'}_{i^{j,k+1}}\left(
M_{i^{j,k+1}}
\right)
\right]
_{s_{j}}.
\]

Let $k\in\bb{\mathbf{c}_{j}}-1$, since $\mathfrak{P}^{(2),i^{j,k},i^{j,k}}$ is, by Proposition~\ref{PDPthSub}, a second-order $c_{i^{j,k}}$-path in $\boldsymbol{\mathcal{A}}^{(2)}$ from $[P_{i^{j,k}}]_{s}$ to $[P_{i^{j,k}+1}]_{s}$, we have that
\begin{itemize}
\item[(i)] $T^{(1)}_{i^{j,k}}\left(
M_{i^{j,k}}
\right)\in [P_{i^{j,k}}]_{s}$ and
\item[(ii)] $T^{(1)}_{i^{j,k}}\left(
N_{i^{j,k}}
\right)\in [P_{i^{j,k}+1}]_{s}$.
\end{itemize}

Let us recall that $T^{(1)}_{i^{j,k}}$, being a first-order translation of type $j$, has the form
\allowdisplaybreaks
\begin{multline*}
T^{(1)}_{i^{j,k}}=\tau^{\mathbf{T}_{\Sigma^{\boldsymbol{\mathcal{A}}}}(X)}\left(
\mathrm{CH}^{(1)}_{s_{0}}\left(\mathfrak{P}_{i^{j,k},0}\right),
\cdots,
\mathrm{CH}^{(1)}_{s_{j-1}}\left(\mathfrak{P}_{i^{j,k},j-1}\right),
\right.
\\
\left.
T^{(1)'}_{i^{j,k}},
\mathrm{CH}^{(1)}_{s_{j+1}}\left(\mathfrak{P}_{i^{j,k},j+1}\right),
\cdots,
\mathrm{CH}^{(1)}_{s_{\bb{\mathbf{s}}-1}}\left(\mathfrak{P}_{i^{j,k},\bb{\mathbf{s}}-1}\right)
\right).
\end{multline*}

Moreover, for the indices after $i^{j,k}$, the corresponding derived first-order translation would not occur at position $j$ until we get to index $i^{j,k+1}$, which would, again, be a first-order translation of type $j$, i.e., a first-order translation of the form
\allowdisplaybreaks
\begin{multline*}
T^{(1)}_{i^{j,k+1}}=\tau^{\mathbf{T}_{\Sigma^{\boldsymbol{\mathcal{A}}}}(X)}\left(
\mathrm{CH}^{(1)}_{s_{0}}\left(\mathfrak{P}_{i^{j,k+1},0}\right),
\cdots,
\mathrm{CH}^{(1)}_{s_{j-1}}\left(\mathfrak{P}_{i^{j,k+1},j-1}\right),
\right.
\\
\left.
T^{(1)'}_{i^{j,k+1}},
\mathrm{CH}^{(1)}_{s_{j+1}}\left(\mathfrak{P}_{i^{j,k+1},j+1}\right),
\cdots,
\mathrm{CH}^{(1)}_{s_{\bb{\mathbf{s}}-1}}\left(\mathfrak{P}_{i^{j,k+1},\bb{\mathbf{s}}-1}\right)
\right).
\end{multline*}

Taking this into account, we would have that, for every index $i\in \bb{\mathbf{c}}$ with $i^{j,k}<i<i^{j,k+1}$, the first-order translation $T^{(1)}_{i}$ would be a first-order translation of type $j_{i}$, with $j_{i}\neq j$, i.e., a first-order translation of the form
\allowdisplaybreaks
\begin{multline*}
T^{(1)}_{i}=\tau^{\mathbf{T}_{\Sigma^{\boldsymbol{\mathcal{A}}}}(X)}\left(
\mathrm{CH}^{(1)}_{s_{0}}\left(\mathfrak{P}_{i,0}\right),
\cdots,
\mathrm{CH}^{(1)}_{s_{j_{i}-1}}\left(\mathfrak{P}_{i,j_{i}-1}\right),
\right.
\\
\left.
T^{(1)'}_{i},
\mathrm{CH}^{(1)}_{s_{j_{i}+1}}\left(\mathfrak{P}_{i,j_{i}+1}\right),
\cdots,
\mathrm{CH}^{(1)}_{s_{\bb{\mathbf{s}}-1}}\left(\mathfrak{P}_{i,\bb{\mathbf{s}}-1}\right)
\right).
\end{multline*}

Since $\mathfrak{P}^{(2)}$ is a coherent second-order path, the following equalities, explicitly located at position $j$, must be fulfilled
\allowdisplaybreaks
\begin{align*}
\left[
T^{(1)'}_{i^{j,k}}\left(
N_{i^{j,k}}
\right)
\right]_{s_{j}}
&=
\left[
\mathrm{CH}^{(1)}_{s_{j}}\left(
\mathfrak{P}_{i^{j,k}+1,j}
\right)
\right]_{s_{j}}
\\&
\qquad\qquad\vdots
\\&=
\left[
\mathrm{CH}^{(1)}_{s_{j}}\left(
\mathfrak{P}_{i^{j,k+1}-1,j}
\right)
\right]_{s_{j}}
\\&=
\left[
T^{(1)'}_{i^{j,k+1}}\left(
M_{i^{j,k+1}}
\right)
\right]_{s_{j}}
\\&=
\left[
P_{j,k+1}
\right]_{s_{j}}.
\end{align*}

This completes the proof.
\end{proof}

We emphasize that the result stated in Lemma~\ref{LDPthExtract} is fundamental for the development of our theory. We will refer to it as the \emph{second-order path extraction algorithm}. To improve the understanding of this lemma we provide the following example showing the computation carried out by the second-order path extraction algorithm in relation to the second-order coherent head-constant stepwise path presented in Example~\ref{EDLoop}.

\begin{example}\label{EDLoopII}[Continuation of Example~\ref{EDLoop}] 
Consider the same second-order rewrite system as in Example~\ref{EHeadCt}.

Consider the second-order path introduced in Example~\ref{EDLoop},
$$
\xymatrix@C=150pt@R=15pt{
\mathfrak{R}^{(2)}: [
\sigma(\mathfrak{q}\circ^{0}\mathfrak{p}, \mathfrak{q}\circ^{0}\mathfrak{p})
]
\ar@2{->}[r]^-{(\mathfrak{p}^{(2)}, 
\sigma(\underline{\quad},z)
\circ^{0}
\sigma(x,\mathfrak{q}\circ^{0}\mathfrak{p}))}
&
[\sigma(\mathfrak{r}, \mathfrak{q}\circ^{0}\mathfrak{p})]
\\
{\qquad\qquad\qquad\qquad\qquad\quad}
\ar@2{->}[r]^-{(\mathfrak{p}^{(2)}, 
\sigma(\mathfrak{r},z) \circ^{0} \sigma(x, \underline{\quad}))}
&
\quad
[\sigma(\mathfrak{r},\mathfrak{r})].
\quad
}
$$
It was stated in Example~\ref{EDLoop} that $\mathfrak{R}^{(2)}$ is a coherent head-constant second-order path associated to the $0$-composition operation symbol $\circ^{0}$.

According to Lemma~\ref{LDPthExtract}, the second-order extraction procedure retrieves the following second-order paths.
$$
\xymatrix@C=100pt{
\mathfrak{R}^{(2)}_{0}: [
\sigma(\mathfrak{q}\circ^{0}\mathfrak{p},z)
]
\ar@2{->}[r]^-{(\mathfrak{p}^{(2)}, 
\sigma(\underline{\quad},z))}
&
[\sigma(\mathfrak{r},z)]
};
$$
$$
\xymatrix@C=100pt{
\mathfrak{R}^{(2)}_{1}: [
\sigma(x,\mathfrak{q}\circ^{0}\mathfrak{p})
]
\ar@2{->}[r]^-{(\mathfrak{p}^{(2)}, 
\sigma(x,\underline{\quad}))}
&
[\sigma(x,\mathfrak{r})]
}.
$$

It can be shown that $\mathfrak{Q}^{(2)}_{0}$ is a second-order path from 
$[\sigma(\mathfrak{q}\circ^{0}\mathfrak{p},z)
]$ to $[\sigma(\mathfrak{r},z)]$ and $\mathfrak{Q}^{(2)}_{1}$ 
a second-order path from 
$[\sigma(x,\mathfrak{q}\circ^{0}\mathfrak{p})
]$ to $[\sigma(x,\mathfrak{r})]$.
\end{example}

The second-order path extraction algorithm from Lemma~\ref{LDPthExtract} can be used to prove different properties on second-order paths following an inductive process. This fact will be used frequently in the following chapters.
\chapter{
\texorpdfstring
{Algebraic structures on $\mathrm{Pth}_{\boldsymbol{\mathcal{A}}^{(2)}}$}
{Algebraic structures on second-order paths}
}\label{S2B}

In this chapter we define a structure of $\Sigma$-algebra on $\mathrm{Pth}_{\boldsymbol{\mathcal{A}}^{(2)}}$, the $S$-sorted set of second-order paths in $\boldsymbol{\mathcal{A}}^{(2)}$, we denote by $\mathbf{Pth}^{(0,2)}_{\boldsymbol{\mathcal{A}}^{(2)}}$ the corresponding  $\Sigma$-algebra. To do this, we introduce for each pair $(\mathbf{s},s)\in S^{\star}\times S$, each operation symbol  $\sigma\in\Sigma_{\mathbf{s},s}$ and each family of second-order paths $(\mathfrak{P}^{(2)}_{j})_{j\in\bb{\mathbf{s}}}$ in $\mathrm{Pth}_{\boldsymbol{\mathcal{A}}^{(2)},\mathbf{s}}$, the second-order path $\sigma^{\mathbf{Pth}_{\boldsymbol{\mathcal{A}}^{(2)}}}((\mathfrak{P}^{(2)}_{j})_{j\in\bb{\mathbf{s}}})$. We next prove that the $(0,1)$-source and the $(0,1)$-target of the $1$-constituents of a second-order path are constant. This allows us to introduce the notions of $(0,2)$-source, $(0,2)$-target and the $(2,0)$-identity second-order paths, denoted by $\mathrm{sc}^{(0,2)}$, $\mathrm{tg}^{(0,2)}$ and $\mathrm{ip}^{(2,0)\sharp}$, respectively. Following this we prove that $(2,0)$-identity second-order paths are completely characterized by their $(0,2)$-source and $(0,2)$-target. We next define a structure of partial $\Sigma^{\boldsymbol{\mathcal{A}}}$-algebra on $\mathrm{Pth}_{\boldsymbol{\mathcal{A}}^{(2)}}$, denoted by $\mathbf{Pth}^{(1,2)}_{\boldsymbol{\mathcal{A}}^{(2)}}$. To do this, we interpret the categorial operations of $\Sigma^{\boldsymbol{\mathcal{A}}}$, including the rewrite rules of $\mathcal{A}$, the $0$-source, the $0$-target and the $0$-composition of second-order paths. We then prove that $\tau^{\mathbf{Pth}_{\boldsymbol{\mathcal{A}}^{(2)}}}((\mathfrak{P}^{(2)}_{j})_{j\in\bb{\mathbf{s}}})$, for $\tau$ being an operation symbol in $\Sigma$ or the $0$-composition is always a coherent head-constant echelonless second-order path and we prove that the second-order path extraction algorithm applied to these compositions always retrieves the original second-order path family $(\mathfrak{P}^{(2)}_{j})_{j\in\bb{\mathbf{s}}}$. Then we prove that the one-step echelonless second-order paths can be completely characterized by means of the extraction algorithm and the compositions introduced above. In addition, we prove that the path $\tau^{\mathbf{Pth}_{\boldsymbol{\mathcal{A}}^{(2)}}}((\mathfrak{P}^{(2)}_{j})_{j\in\bb{\mathbf{s}}})$, for an operation symbol $\tau\in\Sigma^{\boldsymbol{\mathcal{A}}}_{\mathbf{s},s}$, when performed on a family of $(2,[1])$-identity second-order paths $(\mathfrak{P}^{(2)}_{j})_{j\in\bb{\mathbf{s}}}$ is itself a $(2,[1])$-identity second-order path. This allows us to obtain the subalgebra of $(2,[1])$-identity second-order paths, denoted by $\mathrm{ip}^{(2,[1])\sharp}[[\mathbf{PT}_{\boldsymbol{\mathcal{A}}}]]$, of $\mathbf{Pth}^{(1,2)}_{\boldsymbol{\mathcal{A}}^{(2)}}$. Then we show that $\mathrm{ip}^{(2,[1])\sharp}[[\mathbf{PT}_{\boldsymbol{\mathcal{A}}}]]$ belongs to the variety $\mathcal{V}(\boldsymbol{\mathcal{E}}^{\boldsymbol{\mathcal{A}}})$ whose associated category is 
$\mathbf{PAlg}(\boldsymbol{\mathcal{E}}^{\boldsymbol{\mathcal{A}}})$. 
Then we show that the many-sorted mappings $\mathrm{sc}^{([1],2)}$, $\mathrm{tg}^{([1],2)}$ and $\mathrm{ip}^{(2,[1])\sharp}$, of $([1],2)$-source, $([1],2)$-target and $(2,[1])$-identity second-order paths, respectively, are $\Sigma^{\boldsymbol{\mathcal{A}}}$-homomorphisms. In fact, we shown that $\mathrm{ip}^{(2,[1])\sharp}$, the $(2,[1])$-identity second-order path mapping, can be obtained by means of a universal property as the homomorphic extension of the mapping $\mathrm{ip}^{(2,X)}$, that, for every $s\in S$, maps every variable in $X_{s}$ to the $(2,[1])$-identity second-order path on its path term class. 
We then prove that the many-sorted $\Sigma^{\boldsymbol{\mathcal{A}}}$-algebras $\mathbf{Pth}^{(1,2)}_{\boldsymbol{\mathcal{A}}^{(2)}}/{\mathrm{Ker}(\mathrm{sc}^{([1],2)})}$, $\mathbf{Pth}^{(1,2)}_{\boldsymbol{\mathcal{A}}^{(2)}}/{\mathrm{Ker}(\mathrm{tg}^{([1],2)})}$ and $[\mathbf{PT}_{\boldsymbol{\mathcal{A}}}]$ are isomorphic. 

Now, regarding the connection between the layers $2$ and $0$, we prove that the path $\sigma^{\mathbf{Pth}_{\boldsymbol{\mathcal{A}}^{(2)}}}((\mathfrak{P}^{(2)}_{j})_{j\in\bb{\mathbf{s}}})$, for an operation symbol $\sigma\in\Sigma_{\mathbf{s},s}$, when performed on a family of $(2,0)$-identity second-order paths $(\mathfrak{P}^{(2)}_{j})_{j\in\bb{\mathbf{s}}}$ is itself a $(2,0)$-identity second-order path. This allows us to present the subalgebra of $(2,0)$-identity second-order paths, denoted by $\mathrm{ip}^{(2,0)\sharp}[\mathbf{T}_{\Sigma}(X)]$. Then we show that the many-sorted mappings $\mathrm{sc}^{(0,2)}$, $\mathrm{tg}^{(0,2)}$ and $\mathrm{ip}^{(2,0)\sharp}$, of $(0,2)$-source, $(0,2)$-target and $(2,0)$-identity second-order paths, respectively, are $\Sigma$-homomorphisms. In fact, it is shown that $\mathrm{ip}^{(2,0)\sharp}$, the $(2,0)$-identity second-order path mapping, can be obtained by Universal Property as the homomorphic extension of the mapping $\mathrm{ip}^{(2,X)}$, that maps, for every $s\in S$, every variable in $X_{s}$ to the $(2,0)$-identity second-order path on it. We then prove that the many-sorted $\Sigma$-algebras $\mathbf{Pth}^{(0,2)}_{\boldsymbol{\mathcal{A}}^{(2)}}/{\mathrm{Ker}(\mathrm{sc}^{(0,2)})}$, $\mathbf{Pth}^{(0,2)}_{\boldsymbol{\mathcal{A}}^{(2)}}/{\mathrm{Ker}(\mathrm{tg}^{(0,2)})}$ and $\mathbf{T}_{\Sigma}(X)$ are isomorphic.


\section{
\texorpdfstring
{A structure of $\Sigma$-algebra on $\mathrm{Pth}_{\boldsymbol{\mathcal{A}}^{(2)}}$}
{An algebra on second-order paths}
}
We start by defining a structure of $\Sigma$-algebra on the $S$-sorted set $\mathrm{Pth}_{\boldsymbol{\mathcal{A}}^{(2)}}$ of second-order paths. The following proposition mimics Proposition~\ref{PPthAlg}.

\begin{restatable}{proposition}{PDPthAlg}
\label{PDPthAlg} 
\index{path!second-order!$\mathbf{Pth}^{(0,2)}_{\boldsymbol{\mathcal{A}}^{(2)}}$}
The $S$-sorted set $\mathrm{Pth}_{\boldsymbol{\mathcal{A}}^{(2)}}$ is equipped, in a natural way with a structure of $\Sigma$-algebra.
\end{restatable}
\begin{proof}
Let us denote by $\mathbf{Pth}^{(0,2)}_{\boldsymbol{\mathcal{A}}^{(2)}}$ the $\Sigma$-algebra defined as follows:

\textsf{(1)} The underlying $S$-sorted set of $\mathbf{Pth}^{(0,2)}_{\boldsymbol{\mathcal{A}}^{(2)}}$ is $\mathrm{Pth}_{\boldsymbol{\mathcal{A}}^{(2)}}=(\mathrm{Pth}_{\boldsymbol{\mathcal{A}}^{(2)},s})_{s\in S}$.

\textsf{(2)} For every $(\mathbf{s},s)\in S^{\star}\times S$ and every operation symbol $\sigma\in \Sigma_{\mathbf{s},s}$, the operation $\sigma^{\mathbf{Pth}^{(0,2)}_{\boldsymbol{\mathcal{A}}^{(2)}}}$, abbreviated to $\sigma^{\mathbf{Pth}_{\boldsymbol{\mathcal{A}}^{(2)}}}$, from $\mathrm{Pth}_{\boldsymbol{\mathcal{A}}^{(2)},\mathbf{s}}$ to $\mathrm{Pth}_{\boldsymbol{\mathcal{A}}^{(2)},s}$ associated to $\sigma$ assigns to a family of second-order paths $(\mathfrak{P}^{(2)}_{j})_{j\in\bb{\mathbf{s}}}\in\mathrm{Pth}_{\boldsymbol{\mathcal{A}}^{(2)},\mathbf{s}}$, where, for every $j\in\bb{\mathbf{s}}$, $\mathfrak{P}^{(2)}_{j}$ is a $\mathbf{c}_{j}$-second-order path in $\boldsymbol{\mathcal{A}}^{(2)}$ from $[P_{j,0}]_{s_{j}}$ to $[P_{j,\bb{\mathbf{c}_{j}}}]_{s_{j}}$ of the form
$$
\mathfrak{P}^{(2)}_{j}
=
\left(
([P_{j,k}]_{s_{j}})_{k\in\bb{\mathbf{c}_{j}}+1},
(\mathfrak{p}^{(2)}_{j,k})_{k\in\bb{\mathbf{c}_{j}}},
(T^{(1)}_{j,k})_{k\in\bb{\mathbf{c}_{j}}}
\right),
$$
for a unique $\mathbf{c}_{j}$ in $S^{\star}$, and a unique pair of path term classes $[P_{j,0}]_{s_{j}}$ and $[P_{j,\bb{\mathbf{c}_{j}}}]_{s_{j}}$ in $[\mathrm{PT}_{\boldsymbol{\mathcal{A}}}]_{s_{j}}$, precisely the $\mathbf{c}$-second-order path in $\boldsymbol{\mathcal{A}}^{(2)}$ of sort $s$ given by
$$
\sigma^{\mathbf{Pth}_{\boldsymbol{\mathcal{A}}^{(2)}}}
\left(
\left(\mathfrak{P}^{(2)}_{j}\right)_{j\in\bb{\mathbf{s}}}
\right)
=
\left(
([P_{i}]_{s})_{i\in\bb{\mathbf{c}}+1},
(\mathfrak{p}^{(2)}_{i})_{i\in\bb{\mathbf{c}}},
(T^{(1)}_{i})_{i\in\bb{\mathbf{c}}}
\right),
$$
where  $\mathbf{c}=\bigcurlywedge_{j\in\bb{\mathbf{s}}}\mathbf{c}_{j}$
is the concatenation, for $j\in\bb{\mathbf{s}}$, of all words $\mathbf{c}_{j}$.

Before explaining the definition of $\sigma^{\mathbf{Pth}_{\boldsymbol{\mathcal{A}}^{(2)}}}((\mathfrak{P}^{(2)}_{j})_{j\in\bb{\mathbf{s}}})$ let us note the following facts. By construction of $\mathbf{c}$, it happens that $\bb{\mathbf{c}}=\sum_{j\in\bb{\mathbf{s}}}\bb{\mathbf{c}_{j}}$. Hence, the $i$-th letter of $\mathbf{c}$ will be the $k$-th
 letter of the subword $\mathbf{c}_{j}$, for a unique $j\in\bb{\mathbf{s}}$ and a unique $k\in\bb{\mathbf{c}_{j}}$. We will write $i=(j,k)$ to denote this dependency. For the aforementioned indices, we have $c_{i}=c_{j,k}$.
 
 Returning to the definition of $\sigma^{\mathbf{Pth}_{\boldsymbol{\mathcal{A}}^{(2)}}}((\mathfrak{P}^{(2)}_{j})_{j\in\bb{\mathbf{s}}})$, for $i\in\bb{\mathbf{c}}$ with $i=(j,k)$, we define the $1$-constituent at step $i$ of $\sigma^{\mathbf{Pth}_{\boldsymbol{\mathcal{A}}^{(2)}}}((\mathfrak{P}^{(2)}_{j})_{j\in\bb{\mathbf{s}}})$ to be the $\Theta^{[1]}$-class 
\[
[P_{i}]_{s}=
\left[
\sigma^{\mathbf{PT}_{\boldsymbol{\mathcal{A}}}}
\left(
P_{0,\bb{\mathbf{c}_{0}}},
\cdots,
P_{j-1,\bb{\mathbf{c}_{j-1}}},
P_{j,k},
P_{j+1,0}
\cdots,
P_{\bb{\mathbf{s}}-1,0}
\right)\right]_{s}.\]

That is, if $i\in\bb{\mathbf{c}}$ and $i=(j,k)$, then we have that, to the left of position $j$, every subterm is equal to a path term of the last $1$-constituent of the corresponding second-order path, and, to the right of position $j$, every subterm is equal to a path term of the initial $1$-constituent  of the corresponding second-order path.  The $j$-th subterm of $P_{i}$ is the $k$-th term appearing in the second-order path $\mathfrak{P}^{(2)}_{j}$.  Note that, by Definition~\ref{PPTCatAlg}, for every $i\in\bb{\mathbf{c}}$, this is a path term in $\mathrm{PT}_{\boldsymbol{\mathcal{A}}, s}$. Moreover, this path term is well-defined according to Lemma~\ref{LThetaCong} and Definition~\ref{DThetaCong}.  In particular, since $0=(0,0)$, we have that 
\[\left[P_{0}\right]_{\Theta^{[1]}_{2}}=\left[\sigma^{\mathbf{PT}_{\boldsymbol{\mathcal{A}}}}\left(\left(P_{j,0}\right)_{j\in\bb{\mathbf{s}}}\right)\right]_{s}.\]

Finally, for the case of $i=\bb{\mathbf{c}}$, we define
$$
\left[P_{\bb{\mathbf{c}}}
\right]_{s}=\left[
\sigma^{\mathbf{PT}_{\boldsymbol{\mathcal{A}}}}
\left(\left(
P_{j,\bb{\mathbf{c}_{j}}}
\right)_{j\in\bb{\mathbf{s}}}\right)
\right]_{s}.
$$
Note that, in virtue of Proposition~\ref{PPTCatAlg},  this is also a path term class in $[\mathrm{PT}_{\boldsymbol{\mathcal{A}}}]_{s}$. Moreover this path term class is well-defined, according to Lemma~\ref{LThetaCong}. 
 
For $i\in\bb{\mathbf{c}}$ with $i=(j,k)$, we define the second-order rewrite rule $\mathfrak{p}^{(2)}_{i}$ to be equal to $\mathfrak{p}^{(2)}_{j,k}$. That is, the $i$-th rewrite rule of $\sigma^{\mathbf{Pth}_{\boldsymbol{\mathcal{A}}^{(2)}}}((\mathfrak{P}^{(2)}_{j})_{j\in\bb{\mathbf{s}}})$ is equal to the $k$-th second-order rewrite rule of the second-order path $\mathfrak{P}^{(2)}_{j}$. 

Finally, for $i\in\bb{\mathbf{c}}$ with $i=(j,k)$, we define the first-order translation at step $i$ of $\sigma^{\mathbf{Pth}_{\boldsymbol{\mathcal{A}}^{(2)}}}((\mathfrak{P}^{(2)}_{j})_{j\in\bb{\mathbf{s}}})$ to be  equal to
\allowdisplaybreaks
\begin{multline*}
T^{(1)}_{i}=
\sigma^{\mathbf{PT}_{\boldsymbol{\mathcal{A}}}}
\left(
\mathrm{CH}^{(1)}_{s_{0}}\left(
\mathrm{ip}^{(1,X)@}_{s_{0}}\left(
P_{0,\bb{\mathbf{c}_{0}}}
\right)
\right),
\cdots,
\mathrm{CH}^{(1)}_{s_{j-1}}\left(
\mathrm{ip}^{(1,X)@}_{s_{j-1}}\left(
P_{j-1,\bb{\mathbf{c}_{j-1}}}
\right)\right),
\right.
\\
\left.
T^{(1)}_{j,k},
\mathrm{CH}^{(1)}_{s_{j+1}}\left(
\mathrm{ip}^{(1,X)@}_{s_{j+1}}\left(
P_{j+1,0}
\right)\right)
\cdots,
\mathrm{CH}^{(1)}_{s_{\bb{\mathbf{s}}-1}}\left(
\mathrm{ip}^{(1,X)@}_{s_{\bb{\mathbf{s}}-1}}\left(
P_{\bb{\mathbf{s}}-1,0}
\right)\right)
\right).
\end{multline*}

That is, if $i\in\bb{\mathbf{c}}$ and $i=(j,k)$, then we have that, to the left of position $j$, every subterm is equal to the normalized path term of the last $1$-constituent of the corresponding second-order path, and, to the right of position $j$, every subterm is equal to the normalized path term of the initial $1$-constituent  of the corresponding second-order path.  The $j$-th subterm of $T^{(1)}_{i}$ is the $k$-th first-order translation appearing in the second-order path $\mathfrak{P}^{(2)}_{j}$.  Note that, by Definition~\ref{DUTrans}, for every $i\in\bb{\mathbf{c}}$, $T^{(1)}_{i}$ is a first-order translation.

\begin{claim}\label{CDPthSigma}
Let $(\mathbf{s},s)$ be a pair in $S^{\star}\times S$ and, for every $j\in\bb{\mathbf{c}}$, $\mathfrak{P}^{(2)}_{j}$ 
a $\mathbf{c}_{j}$-second-order path in $\mathrm{Pth}_{\boldsymbol{\mathcal{A}}^{(2)},s_{j}}$. Then $\sigma^{\mathbf{Pth}_{\boldsymbol{\mathcal{A}}^{(2)}}}((\mathfrak{P}^{(2)}_{j})_{j\in\bb{\mathbf{s}}})$ is a $\bigcurlywedge_{j\in\bb{\mathbf{s}}}\mathbf{c}_{j}$-path in $\boldsymbol{\mathcal{A}}^{(2)}$ of the form
\begin{multline*}
\sigma^{\mathbf{Pth}_{\boldsymbol{\mathcal{A}}^{(2)}}}
\left(\left(\mathfrak{P}^{(2)}_{j}
\right)_{j\in\bb{\mathbf{s}}}
\right)
\colon
\sigma
^{[\mathbf{PT}_{\boldsymbol{\mathcal{A}}}]}
\left(\left(
\mathrm{sc}^{([1],2)}_{s_{j}}\left(\mathfrak{P}^{(2)}_{j}
\right)
\right)_{j\in\bb{\mathbf{s}}}
\right){\implies}
\\
\sigma
^{[\mathbf{PT}_{\boldsymbol{\mathcal{A}}}]}
\left(\left(
\mathrm{tg}^{([1],2)}_{s_{j}}\left(
\mathfrak{P}^{(2)}_{j}
\right)
\right)_{j\in\bb{\mathbf{s}}}
\right).
\end{multline*}
\end{claim}

Regarding the $([1],2)$-source, the following chain of equalities holds
\allowdisplaybreaks
\begin{align*}
\left[P_{0}
\right]_{s}&
=
\left[\sigma^{\mathbf{PT}_{\boldsymbol{\mathcal{A}}}}
\left(\left(P_{j,0}
\right)_{j\in\bb{\mathbf{s}}}
\right)\right]_{s}
\tag{1}
\\&=
\sigma^{[\mathbf{PT}_{\boldsymbol{\mathcal{A}}}]}
\left(\left(\left[P_{j,0}
\right]_{s_{j}}
\right)_{j\in\bb{\mathbf{s}}}
\right)
\tag{2}
\\&=
\sigma^{[\mathbf{PT}_{\boldsymbol{\mathcal{A}}}]}
\left(\left(\mathrm{sc}_{s_{j}}^{([1],2)}\left(\mathfrak{P}^{(2)}_{j}
\right)\right)_{j\in\bb{\mathbf{s}}}\right).
\tag{3}
\end{align*}

The first equality unpacks the definition of the path term $P_{0}$; the second equality unpacks the interpretation of the operation symbol $\sigma$ in the partial $\Sigma^{\boldsymbol{\mathcal{A}}}$-algebra $[\mathbf{PT}_{\boldsymbol{\mathcal{A}}}]$;  finally, the last equality unpacks the definition of the respective $([1],2)$-sources of the family of second-order paths $(\mathfrak{P}^{(2)}_{j})_{j\in\bb{\mathbf{s}}}$.

Regarding the $([1],2)$-target, the following chain of equalities holds
\allowdisplaybreaks
\begin{align*}
\left[P_{\bb{\mathbf{c}}}
\right]_{s}
&
=
\left[\sigma^{\mathbf{PT}_{\boldsymbol{\mathcal{A}}}}
\left(\left(
P_{j,\bb{\mathbf{c}_{j}}}
\right)_{j\in\bb{\mathbf{s}}}
\right)
\right]_{s}
\tag{1}
\\&=
\sigma
^{[\mathbf{PT}_{\boldsymbol{\mathcal{A}}}]}
\left(\left(
\left[P_{j,\bb{\mathbf{c}_{j}}}
\right]_{s_{j}}
\right)_{j\in\bb{\mathbf{s}}}
\right)
\tag{2}
\\&=
\sigma
^{[\mathbf{PT}_{\boldsymbol{\mathcal{A}}}]}
\left(\left(
\mathrm{tg}^{([1],2)}_{s_{j}}\left(\mathfrak{P}^{(2)}_{j}
\right)\right)_{j\in\bb{\mathbf{s}}}\right).
\tag{3}
\end{align*}

The first equality unpacks the definition of the last path term class in the sequence determined by $\sigma^{\mathbf{Pth}_{\boldsymbol{\mathcal{A}}^{(2)}}}((\mathfrak{P}^{(2)}_{j})_{j\in\bb{\mathbf{s}}})$; the second equality unpacks the interpretation of the operation symbol $\sigma$ in the partial $\Sigma^{\boldsymbol{\mathcal{A}}}$-algebra $[\mathbf{PT}_{\boldsymbol{\mathcal{A}}}]$;  finally, the last equality unpacks the definition of the respective $([1],2)$-targets of the family of second-order paths $(\mathfrak{P}^{(2)}_{j})_{j\in\bb{\mathbf{s}}}$.

Note that for every $i\in\bb{\mathbf{c}}$, there exists $j\in\bb{\mathbf{s}}$ and $k\in\bb{\mathbf{c}_{j}}$ such that $i=(j,k)$ or, equivalently, $c_{i}=c_{j,k}$. Hence, $\mathfrak{p}^{(2)}_{i}=\mathfrak{p}^{(2)}_{j,k}$, where $\mathfrak{p}^{(2)}_{j,k}$ is the $k$-th second-order rewrite rule of the $j$-th second-order path $\mathfrak{P}^{(2)}_{j}$. Let us assume that $\mathfrak{p}^{(2)}_{j,k}=([M_{j,k}]_{c_{j,k}}, [N_{j,k}]_{c_{j,k}})$. Then, since $\mathfrak{P}^{(2)}_{j}$ is a second-order $\mathbf{c}_{j}$-path in $\boldsymbol{\mathcal{A}}^{(2)}$ from $[P_{j,0}]_{s_{j}}$ to $[P_{j,\bb{\mathbf{c}}_{j}}]_{s_{j}}$, it follows that
\begin{multicols}{2}
\begin{itemize}
\item[(i)] $T^{(1)}_{j,k}(M_{j,k})\in [P_{j,k}]_{s_{j}}$;
\item[(ii)] $T^{(1)}_{j,k}(N_{j,k})\in [P_{j,k+1}]_{s_{j}}$.
\end{itemize}
\end{multicols}

Note that the following equality holds
\allowdisplaybreaks
\begin{multline*}
T^{(1)}_{i}\left(
M_{i}
\right)\\=
\sigma^{\mathbf{PT}_{\boldsymbol{\mathcal{A}}}}
\left(
\mathrm{CH}^{(1)}_{s_{0}}\left(
\mathrm{ip}^{(1,X)@}_{s_{0}}\left(
P_{0,\bb{\mathbf{c}_{0}}}
\right)
\right),
\cdots,
\mathrm{CH}^{(1)}_{s_{j-1}}\left(
\mathrm{ip}^{(1,X)@}_{s_{j-1}}\left(
P_{j-1,\bb{\mathbf{c}_{j-1}}}
\right)\right),
\right.
\\
\left.
T^{(1)}_{j,k}\left(
M_{j,k}
\right),
\mathrm{CH}^{(1)}_{s_{j+1}}\left(
\mathrm{ip}^{(1,X)@}_{s_{j+1}}\left(
P_{j+1,0}
\right)\right),
\cdots,
\mathrm{CH}^{(1)}_{s_{\bb{\mathbf{s}}-1}}\left(
\mathrm{ip}^{(1,X)@}_{s_{\bb{\mathbf{s}}-1}}\left(
P_{\bb{\mathbf{s}}-1,0}
\right)\right)
\right).
\end{multline*}

The just stated equality simply unpacks the description of the $i$-th first-order translation $T^{(1)}_{i}$, for an index $i\in \bb{\mathbf{c}}$ with $i=(j,k)$. Note that the following statements hold in virtue of Lemma~\ref{LWCong}
\[
\begin{cases}
\mathrm{CH}^{(1)}_{s_{l}}\left(
\mathrm{ip}^{(1,X)@}_{s_{l}}\left(
P_{l,\bb{\mathbf{c}_{l}}}
\right)\right)
\in \left[P_{l,\bb{\mathbf{c}_{l}}}\right]_{s_{l}}
&\mbox{for every }l\in\bb{\mathbf{s}}\mbox{ with }l<j;
\\
\mathrm{CH}^{(1)}_{s_{l}}\left(
\mathrm{ip}^{(1,X)@}_{s_{l}}\left(
P_{l,0}
\right)\right)
\in \left[P_{l,0}\right]_{s_{l}}
&\mbox{for every }l\in\bb{\mathbf{s}}\mbox{ with }j<l.
\end{cases}
\]

Since $\Theta^{[1]}$ is a $\Sigma^{\boldsymbol{\mathcal{A}}}$-congruence on $\mathbf{PT}_{\boldsymbol{\mathcal{A}}}$ in virtue of Definition~\ref{DThetaCong}, we conclude that 
$
T^{(1)}_{i}\left(
M_{i}
\right)
\in 
\left[
P_{i}
\right]_{s}$.

Now, for the index $i+1$, either (1) the $(i+1)$-th letter of $\mathbf{c}$ is also a letter of the subword $\mathbf{c}_{j}$, and therefore $i+1=(j,k+1)$ or (2) the $(i+1)$-th letter of $\mathbf{c}$ is the first letter of the subword $\mathbf{c}_{j+1}$, hence $i+1=(j+1,0)$ and, necessarily, $i=(j,\bb{\mathbf{c}_{j}}-1)$. These two cases can be handled similarly.

Taking into account the above statements and the fact that, in virtue of Definition~\ref{DThetaCong}, 
$\Theta^{[1]}$ is a $\Sigma^{\boldsymbol{\mathcal{A}}}$-congruence on $\mathbf{PT}_{\boldsymbol{\mathcal{A}}}$, we conclude that 
$
T^{(1)}_{i}\left(
N_{i}
\right)
\in 
\left[
P_{i+1}
\right]_{s}.
$

Thus, $\sigma^{\mathbf{Pth}_{\boldsymbol{\mathcal{A}}^{(2)}}}((\mathfrak{P}^{(2)}_{j})_{j\in\bb{\mathbf{s}}})$ is a second-order $\mathbf{c}$-path in $\boldsymbol{\mathcal{A}}^{(2)}$ from the path term class 
$\sigma
^{[\mathbf{PT}_{\boldsymbol{\mathcal{A}}}]}((
\mathrm{sc}^{([1],2)}_{s_{j}}(\mathfrak{P}^{(2)}_{j})
)_{j\in\bb{\mathbf{s}}})
$ to 
$
\sigma
^{[\mathbf{PT}_{\boldsymbol{\mathcal{A}}}]}((
\mathrm{tg}^{([1],2)}_{s_{j}}(\mathfrak{P}^{(2)}_{j})
)_{j\in\bb{\mathbf{s}}})$.

This finishes the proof.
\end{proof}

\section{
\texorpdfstring
{A structure of partial $\Sigma^{\boldsymbol{\mathcal{A}}}$-algebra on $\mathrm{Pth}_{\boldsymbol{\mathcal{A}}^{(2)}}$}
{A partial algebra on second-order paths}
}

We next equip $\mathrm{Pth}_{\boldsymbol{\mathcal{A}}^{(2)}}$ with a structure of partial $\Sigma^{\boldsymbol{\mathcal{A}}}$-algebra. In this regard, we need to interpret the  categorial operation symbols in 
$\Sigma^{\boldsymbol{\mathcal{A}}}$ as operations on second-order paths. Regarding the constant operation symbols coming from rewrite rules in $\mathcal{A}$, we next introduce the following $S$-sorted mappings.

\begin{restatable}{definition}{DDEchA}
\label{DDEchA}
We will denote by
\begin{enumerate}
\item $\mathrm{ech}^{([1],\mathcal{A})}$ \index{echelon!first-order!$\mathrm{ech}^{([1],\mathcal{A})}$} the $S$-sorted mapping $\mathrm{pr}^{\mathrm{Ker}(\mathrm{CH})}\circ\mathrm{ech}^{(1,\mathcal{A})	}$ from $\mathcal{A}$ to $[\mathrm{Pth}_{\boldsymbol{\mathcal{A}}}]$. Thus, for every sort $s\in S$, $\mathrm{ech}^{([1],\mathcal{A})}$ sends a rewrite rule $\mathfrak{p}\in\mathcal{A}_{s}$ to $[\mathrm{ech}^{(1,\mathcal{A})}(\mathfrak{p})]_{s}$, i.e., the $\mathrm{Ker}(\mathrm{CH}^{(1)})_{s}$-class of the echelon determined by $\mathfrak{p}$.
\item $\eta^{([1],\mathcal{A})}$\index{inclusion!first-order!$\eta^{([1],\mathcal{A})}$} the $S$-sorted mapping $\mathrm{pr}^{\Theta^{[1]}}\circ\eta^{(1,\mathcal{A})}$ from $\mathcal{A}$ to $[\mathrm{PT}_{\boldsymbol{\mathcal{A}}}]$. Thus, for every sort $s\in S$, $\eta^{([1],\mathcal{A})}$ sends a rewrite rule $\mathfrak{p}\in\mathcal{A}_{s}$ to $[\mathfrak{p}^{\mathbf{PT}_{\boldsymbol{\mathcal{A}}}}]_{s}$, i.e., the $\Theta^{[1]}_{s}$-class of the path term determined by $\mathfrak{p}$. 
\item $\mathrm{ech}^{(2,\mathcal{A})}$\index{echelon!second-order!$\mathrm{ech}^{(2,\mathcal{A})}$}  the $S$-sorted mapping from $\mathcal{A}$ to $\mathrm{Pth}_{\boldsymbol{\mathcal{A}}^{(2)}}$ that, for every sort $s\in S$, sends $\mathfrak{p}\in \mathcal{A}_{s}$ to $([\mathfrak{p}^{\mathbf{PT}_{\boldsymbol{\mathcal{A}}}}]_{s},\lambda,\lambda)$, the $(2,[1])$-identity path on $[\mathfrak{p}^{\mathbf{PT}_{\boldsymbol{\mathcal{A}}}}]_{s}$.
\end{enumerate}

The above $S$-sorted mappings are depicted in the diagram of Figure~\ref{FDEchA}.
\end{restatable}

\begin{figure}
\begin{center}
\begin{tikzpicture}
[ACliment/.style={-{To [angle'=45, length=5.75pt, width=4pt, round]}},scale=1.1]
\node[] (p) at (6,-3) [] {$[\mathrm{Pth}_{\boldsymbol{\mathcal{A}}}]$};
\node[] (pt) at (6,-6) [] {$[\mathrm{PT}_{\boldsymbol{\mathcal{A}}}]$};
\node[] (p2) at (6,-9) [] {$\mathrm{Pth}_{\boldsymbol{\mathcal{A}}^{(2)}}$};
\node[] (a) at (0,-3) [] {$\mathcal{A}$};

\draw[ACliment]  (a) to node [above, pos=.7, fill=white]
{$\mathrm{ech}^{([1],\mathcal{A})}$} (p);
\draw[ACliment, bend right=15]  (a) to node  [pos=.6, fill=white] {$\eta^{([1],\mathcal{A})}$} (pt);
\draw[ACliment, bend right=20]  (a) to node [below, pos=.5, fill=white]
{$\mathrm{ech}^{(2,\mathcal{A})}$} (p2);

\node[] (B1) at (6,-4.5)  [] {$\cong$};
\draw[ACliment]  ($(B1)+(-.3,1.2)$) to node [left] {
$\mathrm{CH}^{[1]}$
} ($(B1)+(-.3,-1.2)$);
\draw[ACliment]  ($(B1)+(.3,-1.2)$) to node [right] {
$\mathrm{ip}^{([1],X)@}$
} ($(B1)+(.3,1.2)$);

\node[] (B2) at (6,-7.5)  [] {};
\draw[ACliment]  ($(B2)+(0,1.2)$) to node [above, fill=white] {
$\textstyle {\mathrm{ip}}^{(2,[1])\sharp}$
} ($(B2)+(0,-1.2)$);
\draw[ACliment, bend right]  ($(B2)+(.3,-1.2)$) to node [ below, fill=white] {
$\textstyle {\mathrm{tg}}^{([1],2)}$
} ($(B2)+(.3,1.2)$);
\draw[ACliment, bend left]  ($(B2)+(-.3,-1.2)$) to node [below, fill=white] {
$\textstyle {\mathrm{sc}}^{([1],2)}$
} ($(B2)+(-.3,1.2)$);
\end{tikzpicture}
\end{center}
\caption{Many-sorted mappings relative to $\mathcal{A}$ at layers 1 \& 2.}\label{FDEchA}
\end{figure}

\begin{proposition}\label{PDBasicEqEch} 
For the $S$-sorted mappings defined above, we have that
\begin{multicols}{2}
\begin{itemize}
\item[(i)] $\mathrm{ip}^{(2,[1])\sharp}\circ\eta^{([1],\mathcal{A})}=\mathrm{ech}^{(2,\mathcal{A})};$
\item[(ii)] $\mathrm{sc}^{([1],2)}\circ\mathrm{ech}^{(2,\mathcal{A})}=\eta^{([1],\mathcal{A})};$
\item[(iii)] $\mathrm{tg}^{([1],2)}\circ\mathrm{ech}^{(2,\mathcal{A})}=\eta^{([1],\mathcal{A})};$
\item[(iv)] $\mathrm{CH}^{[1]}\circ\mathrm{ech}^{([1],\mathcal{A})}=\eta^{([1],\mathcal{A})};$
\item[(v)] $\mathrm{ip}^{([1],X)@}\circ\eta^{([1],\mathcal{A})}=\mathrm{ech}^{([1],\mathcal{A})}.$
\end{itemize}
\end{multicols}
\end{proposition}

In order to present the interpretation of the unary operations symbols of $0$-source and $0$-target in $\Sigma^{\boldsymbol{\mathcal{A}}}_{s,s}$, and the binary operation symbols of $0$-composition in $\Sigma^{\boldsymbol{\mathcal{A}}}_{ss,s}$, we introduce the notion of $(0,2)$-source and $(0,2)$-target of  a second-order path. 

In the following proposition we show that every path term class in a second-order path, when interpreted as a path by means of the free completion of the identity mapping, has the same $(0,1)$-source and $(0,1)$-target.

\begin{restatable}{proposition}{PDScTgZ}
\label{PDScTgZ} 
Let $s$ be a sort in $S$, $\mathbf{c}$ a word in $S^{\star}$ and $\mathfrak{P}^{(2)}$ a second-order $\mathbf{c}$-path in $\boldsymbol{\mathcal{A}}^{(2)}$ of the form
$$
\mathfrak{P}^{(2)}
=
\left(
([P_{i}]_{s})_{i\in\bb{\mathbf{c}}+1},
(\mathfrak{p}^{(2)}_{i})_{i\in\bb{\mathbf{c}}},
(T^{(1)}_{i})_{i\in\bb{\mathbf{c}}}
\right).
$$
Then, for every $i,j\in\bb{\mathbf{c}}+1$, we have that
$$
\left(
\mathrm{ip}^{(1,X)@}_{s}
\left(P_{i}
\right),
\mathrm{ip}^{(1,X)@}_{s}\left(
P_{j}
\right)
\right)\in
\mathrm{Ker}\left({\mathrm{sc}}^{(0,1)}\right)_{s}
\cap\mathrm{Ker}\left({\mathrm{tg}}^{(0,1)}\right)_{s}.
$$
\end{restatable}
\begin{proof}
We prove it by induction on $\bb{\mathfrak{P}^{(2)}}$, the length of $\mathfrak{P}^{(2)}$.

\textsf{Base step of the induction.}

If $\bb{\mathfrak{P}^{(2)}}=0$, then $\mathfrak{P}^{(2)}$ is a $(2,[1])$-identity second-order path on the path term class $[P_{0}]_{s}$. In this case, $i=j$ and the statement trivially holds.

\textsf{Inductive step of the induction.}

Assume that the statement holds for second-order paths of length $n\in\mathbb{N}$. That is, if $
\mathfrak{P}^{(2)}
$ is  a second-order $\mathbf{c}$-path in $\boldsymbol{\mathcal{A}}^{(2)}$ of the form
$$
\mathfrak{P}^{(2)}
=
\left(
([P_{i}]_{s})_{i\in\bb{\mathbf{c}}+1},
(\mathfrak{p}^{(2)}_{i})_{i\in\bb{\mathbf{c}}},
(T^{(1)}_{i})_{i\in\bb{\mathbf{c}}}
\right)
$$
and the length of $\bb{\mathfrak{P}^{(2)}}$ is equal to $n$, then for every $i,j\in\bb{\mathbf{c}}+1$, we have that
$$
\left(\mathrm{ip}^{(1,X)@}_{s}\left(
P_{i}
\right),
\mathrm{ip}^{(1,X)@}_{s}\left(
P_{j}
\right)
\right)\in
\mathrm{Ker}\left(
{\mathrm{sc}}^{(0,1)}
\right)_{s}
\cap\mathrm{Ker}\left(
{\mathrm{tg}}^{(0,1)}
\right)_{s}.
$$

Let $
\mathfrak{P}^{(2)}
$ be a second-order $\mathbf{c}$-path in $\boldsymbol{\mathcal{A}}^{(2)}$ of length $n+1$ of the form $$
\mathfrak{P}^{(2)}
=
\left(
([P_{i}]_{s})_{i\in\bb{\mathbf{c}}+1},
(\mathfrak{p}^{(2)}_{i})_{i\in\bb{\mathbf{c}}},
(T^{(1)}_{i})_{i\in\bb{\mathbf{c}}}
\right).
$$
Then the subpath $\mathfrak{P}^{(2),0,\bb{\mathbf{c}}-2}$ is a second-order path from in $\boldsymbol{\mathcal{A}}^{(2)}$ of length $n$. By the inductive hypothesis, for every $i,j\in\bb{\mathbf{c}}$, it is the case that
$$
\left(
\mathrm{ip}^{(1,X)@}_{s}
\left(
P_{i}
\right),
\mathrm{ip}^{(1,X)@}_{s}\left(
P_{j}
\right)
\right)\in
\mathrm{Ker}\left(
{\mathrm{sc}}^{(0,1)}
\right)_{s}
\cap\mathrm{Ker}\left(
{\mathrm{tg}}^{(0,1)}
\right)_{s}.
$$

Now consider the final one-step subpath $\mathfrak{P}^{(2),\bb{\mathbf{c}}-1, \bb{\mathbf{c}}-1}$, that is 
$$
\mathfrak{P}^{(2),\bb{\mathbf{c}}-1, \bb{\mathbf{c}}-1}
\colon
\xymatrix@C=85pt{
\left[P_{\bb{\mathbf{c}}-1}\right]_{s}
\ar@{=>}[r]^-{\text{\Small{$\left(\mathfrak{p}^{(2)}_{\bb{\mathbf{c}}-1}, T^{(1)}_{\bb{\mathbf{c}}-1}\right)$
}}}
&
\left[P_{\bb{\mathbf{c}}}\right]_{s}.
}
$$
Then, if $\mathfrak{p}^{(2)}_{\bb{\mathbf{c}}-1}$, the last second-order rewrite rule occurring in $\mathfrak{P}^{(2)}$, has the form $([M_{\bb{\mathbf{c}}-1}]_{c_{\bb{\mathbf{c}}-1}},[N_{\bb{\mathbf{c}}-1}]_{c_{\bb{\mathbf{c}}-1}})$, then
\begin{multicols}{2}
\begin{itemize}
\item[(i)] $T^{(1)}_{\bb{\mathbf{c}}-1}(
M_{\bb{\mathbf{c}}-1}
)
\in [P_{\bb{\mathbf{c}}-1}]_{s};
$
\item[(ii)] 
$T^{(1)}_{\bb{\mathbf{c}}-1}(
N_{\bb{\mathbf{c}}-1}
)
\in [P_{\bb{\mathbf{c}}}]_{s}.
$
\end{itemize}
\end{multicols}

According to Lemma~\ref{LUTransWD}, we have that the pair
\[
\left(
\mathrm{ip}^{(1,X)@}_{s}\left(
T^{(1)}_{\bb{\mathbf{c}}-1}\left(
M_{\bb{\mathbf{c}}-1}
\right)
\right),
\mathrm{ip}^{(1,X)@}_{s}\left(
T^{(1)}_{\bb{\mathbf{c}}-1}\left(
N_{\bb{\mathbf{c}}-1}
\right)
\right)
\right)
\]
is in $\mathrm{Ker}(
\mathrm{sc}^{(0,1)}
)_{s}\cap\mathrm{Ker}(
\mathrm{tg}^{(0,1)}
)_{s}$.

Moreover, by Corollary~\ref{CPTScTg}, we have that the pair
\[
\left(
\mathrm{ip}^{(1,X)@}_{s}\left(
T^{(1)}_{\bb{\mathbf{c}}-1}\left(
M_{\bb{\mathbf{c}}-1}
\right)
\right),
\mathrm{ip}^{(1,X)@}_{s}\left(
P_{\bb{\mathbf{c}}-1}
\right)
\right);
\]
and the pair
\[
\left(
\mathrm{ip}^{(1,X)@}_{s}\left(
T^{(1)}_{\bb{\mathbf{c}}-1}\left(
N_{\bb{\mathbf{c}}-1}
\right)
\right),
\mathrm{ip}^{(1,X)@}_{s}\left(
P_{\bb{\mathbf{c}}}
\right)
\right),
\]
are also  in $\mathrm{Ker}(
\mathrm{sc}^{(0,1)}
)_{s}\cap\mathrm{Ker}(
\mathrm{tg}^{(0,1)}
)_{s}$.

All in all, we can affirm that 
$$
\left(
\mathrm{ip}^{(1,X)@}_{s}\left(
P_{\bb{\mathbf{c}}-1}
\right),
\mathrm{ip}^{(1,X)@}_{s}\left(
P_{\bb{\mathbf{c}}}\right)\right)
\in\mathrm{Ker}\left(
\mathrm{sc}^{(0,1)}
\right)_{s}\cap\mathrm{Ker}
\left(
\mathrm{tg}^{(0,1)}
\right)_{s}.
$$

Therefore, for every $i,j\in\bb{\mathbf{c}}+1$, we have that 
$$
\left(
\mathrm{ip}^{(1,X)@}_{s}\left(
P_{i}
\right),
\mathrm{ip}^{(1,X)@}_{s}\left(
P_{j}
\right)
\right)\in
\mathrm{Ker}\left(
{\mathrm{sc}}^{(0,1)}
\right)_{s}
\cap\mathrm{Ker}\left(
{\mathrm{tg}}^{(0,1)}
\right)_{s}.
$$

This finishes the proof.
\end{proof}

\begin{corollary}\label{CDScTgZ} Let $s$ be a sort in $S$ and  
$
\mathfrak{P}^{(2)}
$ be a second-order path in $\mathrm{Pth}_{\boldsymbol{\mathcal{A}}^{(2)},s}$. Then 
the following equalities hold
\allowdisplaybreaks
\begin{align*}
\mathrm{sc}^{(0,[1])}_{s}\left(
\mathrm{ip}^{([1],X)@}_{s}\left(
\mathrm{sc}^{([1],2)}_{s}\left(
\mathfrak{P}^{(2)}
\right)\right)\right)&=
\mathrm{sc}^{(0,[1])}_{s}\left(
\mathrm{ip}^{([1],X)@}_{s}\left(
\mathrm{tg}^{([1],2)}_{s}\left(
\mathfrak{P}^{(2)}
\right)\right)\right);
\\
\mathrm{tg}^{(0,[1])}_{s}\left(
\mathrm{ip}^{([1],X)@}_{s}\left(
\mathrm{sc}^{([1],2)}_{s}\left(
\mathfrak{P}^{(2)}
\right)\right)\right)&=
\mathrm{tg}^{(0,[1])}_{s}\left(
\mathrm{ip}^{([1],X)@}_{s}\left(
\mathrm{tg}^{([1],2)}_{s}\left(
\mathfrak{P}^{(2)}
\right)\right)\right).
\end{align*}
\end{corollary}

\begin{corollary}\label{CDScTgZII} Let $s$ be a sort in $S$ and  
$
\mathfrak{P}^{(2)}
$ be a second-order path in $\mathrm{Pth}_{\boldsymbol{\mathcal{A}}^{(2)},s}$. Then the following equalities hold
\begin{align*}
\mathrm{sc}^{0[\mathbf{PT}_{\boldsymbol{\mathcal{A}}}]}_{s}\left(
\mathrm{sc}^{([1],2)}_{s}\left(
\mathfrak{P}^{(2)}
\right)\right)
&=
\mathrm{sc}^{0[\mathbf{PT}_{\boldsymbol{\mathcal{A}}}]}_{s}\left(
\mathrm{tg}^{([1],2)}_{s}\left(
\mathfrak{P}^{(2)}
\right)\right);
\\
\mathrm{tg}^{0[\mathbf{PT}_{\boldsymbol{\mathcal{A}}}]}_{s}\left(
\mathrm{sc}^{([1],2)}_{s}\left(
\mathfrak{P}^{(2)}
\right)\right)
&=
\mathrm{tg}^{0[\mathbf{PT}_{\boldsymbol{\mathcal{A}}}]}_{s}\left(
\mathrm{tg}^{([1],2)}_{s}\left(
\mathfrak{P}^{(2)}
\right)\right).
\end{align*}
\end{corollary}
\begin{proof} 
Let us assume that $\mathfrak{P}^{(2)}$ is a second-order path of the form
$$
\mathfrak{P}^{(2)}
=
\left(
([P_{i}]_{s})_{i\in\bb{\mathbf{c}}+1},
(\mathfrak{p}^{(2)}_{i})_{i\in\bb{\mathbf{c}}},
(T^{(1)}_{i})_{i\in\bb{\mathbf{c}}}
\right).
$$

Note that the following chain of equalities hold
\allowdisplaybreaks
\begin{align*}
\mathrm{sc}^{0[\mathbf{PT}_{\boldsymbol{\mathcal{A}}}]}_{s}\left(
\mathrm{sc}^{([1],2)}_{s}\left(
\mathfrak{P}^{(2)}
\right)\right) &=
\mathrm{sc}^{0[\mathbf{PT}_{\boldsymbol{\mathcal{A}}}]}_{s}\left(
[P_{0}]_{s}
\right)
\tag{1}
\\&=
\left[
\mathrm{sc}^{0\mathbf{PT}_{\boldsymbol{\mathcal{A}}}}_{s}
\left(
P_{0}
\right)
\right]_{s}
\tag{2}
\\&=
\left[
\mathrm{CH}^{(1)}_{s}\left(
\mathrm{ip}^{(1,X)@}_{s}\left(
\mathrm{sc}^{0\mathbf{PT}_{\boldsymbol{\mathcal{A}}}}_{s}\left(
P_{0}
\right)\right)\right)
\right]_{s}
\tag{3}
\\&=
\left[
\mathrm{CH}^{(1)}_{s}\left(
\mathrm{sc}^{0\mathbf{Pth}_{\boldsymbol{\mathcal{A}}}}_{s}\left(
\mathrm{ip}^{(1,X)@}_{s}\left(
P_{0}
\right)\right)\right)
\right]_{s}
\tag{4}
\\&=
\left[
\mathrm{CH}^{(1)}_{s}\left(
\mathrm{ip}^{(1,0)\sharp}_{s}\left(
\mathrm{sc}^{(0,1)}_{s}\left(
\mathrm{ip}^{(1,X)@}_{s}\left(
P_{0}
\right)\right)\right)\right)
\right]_{s}
\tag{5}
\\&=
\left[
\mathrm{CH}^{(1)}_{s}\left(
\mathrm{ip}^{(1,0)\sharp}_{s}\left(
\mathrm{sc}^{(0,1)}_{s}\left(
\mathrm{ip}^{(1,X)@}_{s}\left(
P_{\bb{\mathbf{c}}}
\right)\right)\right)\right)
\right]_{s}
\tag{6}
\\&=
\left[
\mathrm{CH}^{(1)}_{s}\left(
\mathrm{sc}^{0\mathbf{Pth}_{\boldsymbol{\mathcal{A}}}}_{s}\left(
\mathrm{ip}^{(1,X)@}_{s}\left(
P_{\bb{\mathbf{c}}}
\right)\right)\right)
\right]_{s}
\tag{7}
\\&=
\left[
\mathrm{CH}^{(1)}_{s}\left(
\mathrm{ip}^{(1,X)@}_{s}\left(
\mathrm{sc}^{0\mathbf{PT}_{\boldsymbol{\mathcal{A}}}}_{s}\left(
P_{\bb{\mathbf{c}}}
\right)\right)\right)
\right]_{s}
\tag{8}
\\&=
\left[
\mathrm{sc}^{0\mathbf{PT}_{\boldsymbol{\mathcal{A}}}}_{s}
\left(
P_{\bb{\mathbf{c}}}
\right)
\right]_{s}
\tag{9}
\\&=
\mathrm{sc}^{0[\mathbf{PT}_{\boldsymbol{\mathcal{A}}}]}_{s}\left(
[P_{\bb{\mathbf{c}}}]_{s}
\right)
\tag{10}
\\&=
\mathrm{sc}^{0[\mathbf{PT}_{\boldsymbol{\mathcal{A}}}]}_{s}\left(
\mathrm{tg}^{([1],2)}_{s}\left(
\mathfrak{P}^{(2)}
\right)\right).
\tag{11}
\end{align*}

In the just stated chain of equalities, the first equality unpacks the $([1],2)$-source of $\mathfrak{P}^{(2)}$; the second equality applies the interpretation of the $0$-source operation symbol in the many-sorted partial $\Sigma^{\boldsymbol{\mathcal{A}}}$-algebra $[\mathbf{PT}_{\boldsymbol{\mathcal{A}}}]$ according to Proposition~\ref{PPTQCatAlg}; the third equality follows from Lemma~\ref{LWCong}; the fourth equality follows from the fact that $\mathrm{ip}^{(1,X)@}$ is a $\Sigma^{\boldsymbol{\mathcal{A}}}$-homomorphism according to Definition~\ref{DIp}; the fifth equality unpacks the description of the $0$-source operation symbol in the many-sorted partial $\Sigma^{\boldsymbol{\mathcal{A}}}$-algebra $\mathbf{Pth}_{\boldsymbol{\mathcal{A}}}$ according to Definition~\ref{PPthCatAlg}; the sixth equality follows from Proposition~\ref{PDScTgZ}; the seventh equality recovers the description of the $0$-source operation symbol in the many-sorted partial $\Sigma^{\boldsymbol{\mathcal{A}}}$-algebra $\mathbf{Pth}_{\boldsymbol{\mathcal{A}}}$ according to Definition~\ref{PPthCatAlg}; the eighth equality follows from the fact that $\mathrm{ip}^{(1,X)@}$ is a $\Sigma^{\boldsymbol{\mathcal{A}}}$-homomorphism according to Definition~\ref{DIp}; the ninth equality follows from Lemma~\ref{LWCong}; the tenth equality recovers  the interpretation of the $0$-source operation symbol in the many-sorted partial $\Sigma^{\boldsymbol{\mathcal{A}}}$-algebra $[\mathbf{PT}_{\boldsymbol{\mathcal{A}}}]$ according to Proposition~\ref{PPTQCatAlg}; finally, the last equality recovers  the $([1],2)$-target of $\mathfrak{P}^{(2)}$.

The second statement can be proved in a similar manner.

This completes the proof.
\end{proof}

We are now in position to define the $(0,2)$-source and the $(0,2)$-target of a second-order path.
We will also introduce the $(2,0)$-identity second-order path on a term.

\begin{restatable}{definition}{DDScTgZ}
\label{DDScTgZ} Let $s$ be a sort in $S$ and  
$
\mathfrak{P}^{(2)}
$  a second-order path in  $\mathrm{Pth}_{\boldsymbol{\mathcal{A}}^{(2)},s}$. Then 
the terms 
$\mathrm{sc}^{(0,[1])}_{s}(\mathrm{ip}^{([1],X)@}_{s}(\mathrm{sc}^{([1],2)}_{s}(\mathfrak{P}^{(2)})))$ and 
$\mathrm{tg}^{(0,[1])}_{s}(\mathrm{ip}^{([1],X)@}_{s}(\mathrm{tg}^{([1],2)}_{s}(\mathfrak{P}^{(2)})))$ 
in $\mathrm{T}_{\Sigma}(X)_{s}$ will be called the \emph{$(0,2)$-source} and the \emph{$(0,2)$-target} of the second-order path $\mathfrak{P}^{(2)}$, respectively. 

We will denote by
\begin{enumerate}
\item $\mathrm{sc}^{(0,2)}$\index{source!second-order!$\mathrm{sc}^{(0,2)}$} the $S$-sorted mapping $\mathrm{sc}^{(0,[1])}\circ\mathrm{ip}^{([1],X)@}\circ\mathrm{sc}^{([1],2)}$ from $\mathrm{Pth}_{\boldsymbol{\mathcal{A}}^{(2)}}$ to $\mathrm{T}_{\Sigma}(X)$. Thus $\mathrm{sc}^{(0,2)}$ sends a second-order path to its $(0,2)$-source.
\item $\mathrm{tg}^{(0,2)}$\index{target!second-order!$\mathrm{tg}^{(0,2)}$}  the $S$-sorted mapping $\mathrm{tg}^{(0,[1])}\circ\mathrm{ip}^{([1],X)@}\circ\mathrm{tg}^{([1],2)}$ from $\mathrm{Pth}_{\boldsymbol{\mathcal{A}}^{(2)}}$ to $\mathrm{T}_{\Sigma}(X)$. Thus $\mathrm{tg}^{(0,2)}$ sends a second-order path to its $(0,2)$-target.
\item $\mathrm{ip}^{(2,0)\sharp}$\index{identity!second-order!$\mathrm{ip}^{(2,0)\sharp}$}  the $S$-sorted mapping $\mathrm{ip}^{(2,[1])\sharp}\circ\mathrm{CH}^{[1]}\circ\mathrm{ip}^{([1],0)\sharp}$ from $\mathrm{T}_{\Sigma}(X)$ to $\mathrm{Pth}_{\boldsymbol{\mathcal{A}}^{(2)}}$. Thus $\mathrm{ip}^{(2,0)\sharp}$ sends a term to its $(2,0)$-identity second-order path.
\end{enumerate} 

\end{restatable}

The aforementioned $S$-sorted mappings are depicted in the diagram of Figure~\ref{FDScTgIpZ}. 

\begin{center}
\begin{figure}
\begin{tikzpicture}
[ACliment/.style={-{To [angle'=45, length=5.75pt, width=4pt, round]}},scale=1.1]
\node[] (xoq) at (2,0) [] {$X$};
\node[] (txoq) at (6,0) [] {$\mathrm{T}_{\Sigma}(X)$};
\node[] (tx2) at (9,0) [] {$\mathrm{T}_{\Sigma}(X)$};
\node[] (p) at (6,-3) [] {$[\mathrm{Pth}_{\boldsymbol{\mathcal{A}}}]$};
\node[] (pt) at (6,-6) [] {$[\mathrm{PT}_{\boldsymbol{\mathcal{A}}}]$};
\node[] (p2) at (6,-9) [] {$\mathrm{Pth}_{\boldsymbol{\mathcal{A}}^{(2)}}$};
\node[] (p22) at (9,-9) [] {$\mathrm{Pth}_{\boldsymbol{\mathcal{A}}^{(2)}}$};

\draw[-,double, double equal sign distance] (txoq) to node [] {} (tx2);
\draw[-,double, double equal sign distance] (p2) to node [] {} (p22);

\draw[ACliment]  (xoq) to node [above]
{$\eta^{(0,X)}$} (txoq);
\draw[ACliment, bend right=10]  (xoq) to node  [midway, fill=white] {$\mathrm{ip}^{([1],X)}$} (p);
\draw[ACliment, bend right=15]  (xoq) to node [midway, fill=white]
{$\eta^{([1],X)}$} (pt);
\draw[ACliment, bend right=20]  (xoq) to node  [below left] {$\mathrm{ip}^{(2,X)}$} (p2);

\node[] (B0) at (6,-1.5)  [] {};
\draw[ACliment]  ($(B0)+(0,1.2)$) to node [above, fill=white] {
$\textstyle \mathrm{ip}^{([1],0)\sharp}$
} ($(B0)+(0,-1.2)$);
\draw[ACliment, bend right]  ($(B0)+(.3,-1.2)$) to node [ below, fill=white] {
$\textstyle \mathrm{tg}^{(0,[1])}$
} ($(B0)+(.3,1.2)$);
\draw[ACliment, bend left]  ($(B0)+(-.3,-1.2)$) to node [below, fill=white] {
$\textstyle \mathrm{sc}^{(0,[1])}$
} ($(B0)+(-.3,1.2)$);

\node[] (B1) at (6,-4.5)  [] {$\cong$};
\draw[ACliment]  ($(B1)+(-.3,1.2)$) to node [left] {
$\textstyle \mathrm{CH}^{[1]}$
} ($(B1)+(-.3,-1.2)$);
\draw[ACliment]  ($(B1)+(.3,-1.2)$) to node [right] {
$\textstyle\mathrm{ip}^{([1],X)@}$
} ($(B1)+(.3,1.2)$);

\node[] (B2) at (6,-7.5)  [] {};
\draw[ACliment]  ($(B2)+(0,1.2)$) to node [above, fill=white] {
$\textstyle {\mathrm{ip}}^{(2,[1])\sharp}$
} ($(B2)+(0,-1.2)$);
\draw[ACliment, bend right]  ($(B2)+(.3,-1.2)$) to node [ below, fill=white] {
$\textstyle {\mathrm{tg}}^{([1],2)}$
} ($(B2)+(.3,1.2)$);
\draw[ACliment, bend left]  ($(B2)+(-.3,-1.2)$) to node [below, fill=white] {
$\textstyle {\mathrm{sc}}^{([1],2)}$
} ($(B2)+(-.3,1.2)$);

\node[] (B3) at (9,-4.5)  [] {};
\draw[ACliment]  ($(B3)+(0,4.2)$) to node [above, fill=white] {
$\textstyle {\mathrm{ip}}^{(2,0)\sharp}$
} ($(B3)+(0,-4.2)$);
\draw[ACliment, bend right=10]  ($(B3)+(.3,-4.2)$) to node [ below, fill=white] {
$\textstyle {\mathrm{tg}}^{(0,2)}$
} ($(B3)+(.3,4.2)$);
\draw[ACliment, bend left=10]  ($(B3)+(-.3,-4.2)$) to node [below, fill=white] {
$\textstyle {\mathrm{sc}}^{(0,2)}$
} ($(B3)+(-.3,4.2)$);
\end{tikzpicture}
\caption{Many-sorted mappings relative to $X$ at layers 0, 1 \& 2.}\label{FDScTgIpZ}
\end{figure}
\end{center}

On the basis of the definition of the previous mappings, we now introduce the notion of $(2,0)$-identity second-order paths.
\begin{definition}\label{DDIpZ} 
Let $s$ be a sort in $S$ and let $\mathfrak{P}^{(2)}$ be a second-order path in $\mathrm{Pth}_{\boldsymbol{\mathcal{A}}^{(2)},s}$. We will say that $\mathfrak{P}^{(2)}$ is a $(2,0)$-\emph{identity second-order path} if $\mathfrak{P}^{(2)}=\mathrm{ip}^{(2,0)\sharp}_{s}(P)$, for some term $P\in\mathrm{T}_{\Sigma}(X)_{s}$.
\end{definition}

In the following proposition we state that the $(2,0)$-identity second-order paths are essentially $([1],0)$-identity paths when these are lifted to the second level as $(2,[1])$-identity second-order paths.

\begin{restatable}{proposition}{PDIpZ}
\label{PDIpZ} Let $s$ be a sort in $S$ and let $\mathfrak{P}^{(2)}$ be a second-order path in $\mathrm{Pth}_{\boldsymbol{\mathcal{A}}^{(2)},s}$, then the following statements are equivalent
\begin{enumerate}
\item $\mathfrak{P}^{(2)}$ is a $(2,0)$-identity second-order path;
\item The following two conditions hold
\begin{enumerate}
\item[(2.1)] $\mathfrak{P}^{(2)}$ is a $(2,[1])$-identity second-order path;
\item[(2.2)] $\mathrm{ip}^{([1],X)@}_{s}(\mathrm{sc}^{([1],2)}_{s}(\mathfrak{P}^{(2)}))$ is a $([1],0)$-identity path.
\end{enumerate}
\end{enumerate} 
\end{restatable}

\begin{proof}
Let us assume that $\mathfrak{P}^{(2)}$ is a $(2,0)$-identity second-order path. Then, by Definition~\ref{DDIpZ}, there exists a term $P\in\mathrm{T}_{\Sigma}(X)_{s}$ for which $\mathfrak{P}^{(2)}=\mathrm{ip}^{(2,0)\sharp}_{s}(P)$.

Taking into account Definition~\ref{DDScTgZ}, we have that 
\[
\mathfrak{P}^{(2)}=
\mathrm{ip}^{(2,0)\sharp}_{s}(P)=
\mathrm{ip}^{(2,[1])\sharp}_{s}\left(
\mathrm{CH}^{[1]}_{s}\left(
\mathrm{ip}^{([1],0)\sharp}_{s}\left(
P
\right)\right)\right).
\]

From the above equality, we conclude that $\mathfrak{P}^{(2)}$ is a $(2,[1])$-identity second-order path. Moreover, 
\[
\mathrm{sc}^{([1],2)}_{s}\left(
\mathfrak{P}^{(2)}
\right)=
\mathrm{CH}^{[1]}_{s}\left(
\mathrm{ip}^{([1],0)\sharp}_{s}\left(
P
\right)\right)
\]

Therefore, the following chain of equalities hold
\begin{align*}
\mathrm{ip}^{([1],X)@}_{s}\left(
\mathrm{sc}^{([1],2)}_{s}\left(
\mathfrak{P}^{(2)}
\right)
\right)
&=
\mathrm{ip}^{([1],X)@}_{s}\left(
\mathrm{CH}^{[1]}_{s}\left(
\mathrm{ip}^{([1],0)\sharp}_{s}\left(
P
\right)\right)
\right)
\tag{1}
\\&=
\mathrm{ip}^{([1],0)\sharp}_{s}\left(
P
\right).
\tag{2}
\end{align*}

In the just stated chain of equalities, the first equality follows from the description of the $([1],2)$-source of $\mathfrak{P}^{(2)}$, and the second equality follows from Theorem~\ref{TIso}.

Therefore, $\mathrm{ip}^{([1],X)@}_{s}(\mathrm{sc}^{([1],2)}_{s}(\mathfrak{P}^{(2)}))$ is a $([1],0)$-identity path.

Now, for the other implication, assume that $\mathfrak{P}^{(2)}$ is a $(2,[1])$-identity second-order path and $\mathrm{ip}^{([1],X)@}_{s}(\mathrm{sc}^{([1],2)}_{s}(\mathfrak{P}^{(2)}))$ is a $([1],0)$-identity path. Since $\mathfrak{P}^{(2)}$ is a $(2,[1])$-identity second-order path, it can be written as 
\[
\mathfrak{P}^{(2)}=\mathrm{ip}^{(2,[1])\sharp}_{s}\left([P]_{s}\right),
\]
for some path term $P$ in $\mathrm{PT}_{\boldsymbol{\mathcal{A}}, s}$. Now, since $\mathrm{ip}^{([1],X)@}_{s}(\mathrm{sc}^{([1],2)\sharp}_{s}(\mathfrak{P}^{(2)}))$ is a $([1],0)$-identity path, there exists a term $Q\in \mathrm{T}_{\Sigma}(X)_{s}$ for which 
\[
\mathrm{ip}^{([1],X)@}_{s}\left(
\mathrm{sc}^{([1],2)\sharp}_{s}\left(\mathfrak{P}^{(2)}\right)
\right)=\mathrm{ip}^{([1],0)\sharp}_{s}(Q).
\]

In this regard, note that $\mathrm{sc}^{([1],2)\sharp}_{s}(\mathfrak{P}^{(2)})=[P]_{s}$. Therefore, we have that 
\[
\mathrm{ip}^{([1],X)@}_{s}\left(
[P]_{s}
\right)=\mathrm{ip}^{([1],0)\sharp}_{s}(Q).
\]

Taking this into account, the following chain of equalities holds
\allowdisplaybreaks
\begin{align*}
[P]_{s}&=\mathrm{CH}^{[1]}_{s}\left(
\mathrm{ip}^{([1],X)@}_{s}\left(
[P]_{s}
\right)\right)
\tag{1}
\\&=
\mathrm{CH}^{[1]}_{s}\left(
\mathrm{ip}^{([1],0)\sharp}_{s}(Q)
\right).
\tag{2}
\end{align*}

In the just stated chain of equalities, the first equality follows from Theorem~\ref{TIso}, and the second equality follows from the discussion above.

All in all, we conclude that
\[
\mathfrak{P}^{(2)}=\mathrm{ip}^{(2,[1])\sharp}_{s}\left(
[P]_{s}
\right)
=
\mathrm{ip}^{(2,[1])\sharp}_{s}\left(
\mathrm{CH}^{[1]}_{s}\left(
\mathrm{ip}^{([1],0)\sharp}_{s}(Q)
\right)
\right)
=
\mathrm{ip}^{(2,0)\sharp}_{s}(Q).
\]

Therefore, $\mathfrak{P}^{(2)}$ is a $(2,0)$-identity second-order path.

This concludes the proof.
\end{proof}

\begin{proposition}\label{PDBasicEqZ} For the $S$-sorted mappings introduced before, we have that
\begin{multicols}{2}
\begin{itemize}
\item[(i)] $\mathrm{ip}^{(2,0)\sharp}\circ\eta^{(0,X)}=\mathrm{ip}^{(2,X)};$
\item[(ii)] $\mathrm{sc}^{(0,2)}\circ\mathrm{ip}^{(2,X)}=\eta^{(0,X)};$
\item[(iii)] $\mathrm{tg}^{(0,2)}\circ\mathrm{ip}^{(2,X)}=\eta^{(0,X)};$
\item[(iv)] $\mathrm{sc}^{(0,2)}\circ\mathrm{ip}^{(2,0)\sharp}=\mathrm{id}^{\mathrm{T}_{\Sigma}(X)};$
\item[(v)] $\mathrm{tg}^{(0,2)}\circ\mathrm{ip}^{(2,0)\sharp}=\mathrm{id}^{\mathrm{T}_{\Sigma}(X)}.$
\end{itemize}
\end{multicols}
\end{proposition}

The diagram in Figure~\ref{FDPth} represents an arbitrary second-order path. The representation is simplified since we are taking, by abuse of notation, the twofold interpretation of path classes in $[\mathrm{Pth}_{\boldsymbol{\mathcal{A}}}]$ and path term classes in $[\mathrm{PT}_{\boldsymbol{\mathcal{A}}}]$, by means of the Curry-Howard isomorphism.

\begin{center}
\begin{figure}
\begin{tikzpicture}
[ACliment/.style={-{To [angle'=45, length=5.75pt, width=4pt, round]}}
, scale=0.8, 
AClimentD/.style={double equal sign distance,
-implies
}
]

\node[] (0) at (0,0) [] {$\mathrm{sc}^{(0,2)}_{s}(\mathfrak{P}^{(2)})$};
\node[] (1) at (5,0) [] {$\mathrm{tg}^{(0,2)}_{s}(\mathfrak{P}^{(2)})$};
\node[] (a) at (2.5,.7) [] {};
\node[] (b) at (2.5,-.7) [] {};

\draw[ACliment, bend left] (0) to node [above] {$\mathrm{sc}^{([1],2)}_{s}(\mathfrak{P}^{(2)})$} (1); 
\draw[ACliment, bend right] (0) to node [below] {$ \mathrm{tg}^{([1],2)}_{s}(\mathfrak{P}^{(2)})$} (1); 
\draw[AClimentD] (a) to node [right] {$ \mathfrak{P}^{(2)}$} (b); 
\end{tikzpicture}
\caption{A simplified representation of a second-order path.}\label{FDPth}
\end{figure}
\end{center}

Taking into account the definition of the $(0,2)$-source and $(0,2)$-target of a second-order path we can introduce a structure of partial $\Sigma^{\boldsymbol{\mathcal{A}}}$-algebra on $\mathrm{Pth}_{\boldsymbol{\mathcal{A}}^{(2)}}$.  Let us recall that the operations associated to operation symbols $\sigma$ coming from the original signature $\Sigma$ follow a leftmost innermost strategy, by analogy with the definitions on Chapter~\ref{S1B}. However, regarding the operation of $0$-composition we reverse the strategy to mimic the usual rightmost writing of the standard $0$-composition.

\begin{restatable}{proposition}{PDPthCatAlg}\label{PDPthCatAlg} 
\index{path!second-order!$\mathbf{Pth}^{(1,2)}_{\boldsymbol{\mathcal{A}}^{(2)}}$}
The $S$-sorted set $\mathrm{Pth}_{\boldsymbol{\mathcal{A}}^{(2)}}$ is equipped, in a natural way with a structure of partial $\Sigma^{\boldsymbol{\mathcal{A}}}$-algebra.
\end{restatable}

\begin{proof}
Let us denote by $\mathbf{Pth}^{(1,2)}_{\boldsymbol{\mathcal{A}}^{(2)}}$ the many-sorted partial $\Sigma^{\boldsymbol{\mathcal{A}}}$-algebra defined as follows:

\textsf{(1)} The underlying $S$-sorted set of $\mathbf{Pth}^{(1,2)}_{\boldsymbol{\mathcal{A}}^{(2)}}$ is $\mathrm{Pth}_{\boldsymbol{\mathcal{A}}^{(2)}}=(\mathrm{Pth}_{\boldsymbol{\mathcal{A}}^{(2)},s})_{s\in S}$.

\textsf{(2)}  For every $(\mathbf{s},s)\in S^{\star}\times S$ and every operation symbol 
$\sigma\in \Sigma_{\mathbf{s},s}$,
 the operation $\sigma^{\mathbf{Pth}^{(1,2)}_{\boldsymbol{\mathcal{A}}^{(2)}}}$ 
is equal to $\sigma^{\mathbf{Pth}^{(0,2)}_{\boldsymbol{\mathcal{A}}^{(2)}}}$, abbreviated to $\sigma^{\mathbf{Pth}_{\boldsymbol{\mathcal{A}}^{(2)}}}$,
i.e., to the interpretation of $\sigma$ in $\mathbf{Pth}^{(0,2)}_{\boldsymbol{\mathcal{A}}^{(2)}}$ 
that, we recall, was stated in Proposition~\ref{PDPthAlg}, where, in addition, we proved that $\sigma^{\mathbf{Pth}_{\boldsymbol{\mathcal{A}}^{(2)}}}$ is well-defined.

\textsf{(3)} For every $s\in S$ and every $\mathfrak{p}\in \mathcal{A}_{s}$, the constant $\mathfrak{p}^{\mathbf{Pth}^{(1,2)}_{\boldsymbol{\mathcal{A}}^{(2)}}}$, abbreviated to  $\mathfrak{p}^{\mathbf{Pth}_{\boldsymbol{\mathcal{A}}^{(2)}}}$, is given by
$$
\mathfrak{p}^{\mathbf{Pth}_{\boldsymbol{\mathcal{A}}^{(2)}}}
=
\mathrm{ech}^{(2,\mathcal{A})}_{s}(\mathfrak{p}),
$$
i.e., by the $(2,[1])$-identity second-order path on $[\mathfrak{p}^{\mathbf{PT}_{\boldsymbol{\mathcal{A}}}}]_{s}$.

\textsf{(4)} For every $s\in S$, the operation of $0$-source $\mathrm{sc}_{s}^{0\mathbf{Pth}^{(1,2)}_{\boldsymbol{\mathcal{A}}^{(2)}}}$, abbreviated to $\mathrm{sc}_{s}^{0\mathbf{Pth}_{\boldsymbol{\mathcal{A}}^{(2)}}}$,  from $\mathrm{Pth}_{\boldsymbol{\mathcal{A}}^{(2)},s}$ to $\mathrm{Pth}_{\boldsymbol{\mathcal{A}}^{(2)},s}$, associated to $\mathrm{sc}^{0}_{s}$ assigns to a second-order path $\mathfrak{P}^{(2)}\in\mathrm{Pth}_{\boldsymbol{\mathcal{A}}^{(2)},s}$
the second-order $\mathbf{c}$-path  given by
\[
\mathrm{sc}_{s}^{0\mathbf{Pth}_{\boldsymbol{\mathcal{A}}^{(2)}}}\left(
\mathfrak{P}^{(2)}\right)=
\mathrm{ip}^{(2,0)\sharp}_{s}\left(
\mathrm{sc}^{(0,2)}_{s}\left(\mathfrak{P}^{(2)}
\right)\right).
\]

\begin{claim}\label{CDPthCatAlgScZ} Let $s$ be a sort in $S$, $\mathbf{c}$ a word in $S^{\star}$, and  $\mathfrak{P}^{(2)}$ a second-order $\mathbf{c}$-path in $\mathrm{Pth}_{\boldsymbol{\mathcal{A}}^{(2)},s}$, then 
$\mathrm{sc}^{0
\mathbf{Pth}_{\boldsymbol{\mathcal{A}}^{(2)}}
}_{s}(\mathfrak{P}^{(2)})$ is the $(2,[1])$-identity second-order path in $\mathrm{Pth}_{\boldsymbol{\mathcal{A}}^{(2)},s}$  given by
\allowdisplaybreaks
\begin{align*}
\mathrm{sc}^{0
\mathbf{Pth}_{\boldsymbol{\mathcal{A}}^{(2)}}
}_{s}\left(\mathfrak{P}^{(2)}\right)
&=
\mathrm{ip}^{(2,[1])\sharp}_{s}\left(
\mathrm{sc}^{0[\mathbf{PT}_{\boldsymbol{\mathcal{A}}}]}_{s}\left(
\mathrm{sc}^{([1],2)}_{s}\left(
\mathfrak{P}^{(2)}
\right)\right)\right)
\\&=
\mathrm{ip}^{(2,[1])\sharp}_{s}\left(
\mathrm{sc}^{0[\mathbf{PT}_{\boldsymbol{\mathcal{A}}}]}_{s}\left(
\mathrm{tg}^{([1],2)}_{s}\left(
\mathfrak{P}^{(2)}
\right)\right)\right).
\end{align*}
\end{claim}

Note that the following chain of equalities holds
\begin{flushleft}
$\mathrm{sc}^{0
\mathbf{Pth}_{\boldsymbol{\mathcal{A}}^{(2)}}
}_{s}\left(\mathfrak{P}^{(2)}\right)$
\allowdisplaybreaks
\begin{align*}
\qquad&=
\mathrm{ip}^{(2,0)\sharp}_{s}\left(
\mathrm{sc}^{(0,2)}_{s}\left(\mathfrak{P}^{(2)}
\right)\right)
\tag{1}
\\
&=
\mathrm{ip}^{(2,[1])\sharp}_{s}\left(
\mathrm{CH}^{[1]}_{s}\left(
\mathrm{ip}^{(1,0)\sharp}_{s}\left(
\mathrm{sc}^{(0,[1])}_{s}\left(
\mathrm{ip}^{([1],X)@}_{s}\left(
\mathrm{sc}^{([1],2)}_{s}\left(
\mathfrak{P}^{(2)}
\right)\right)\right)\right)\right)\right)
\tag{2}
\\
&=
\mathrm{ip}^{(2,[1])\sharp}_{s}\left(
\mathrm{CH}^{[1]}_{s}\left(
\mathrm{sc}^{0[\mathbf{Pth}_{\boldsymbol{\mathcal{A}}}]}_{s}\left(
\mathrm{ip}^{([1],X)@}_{s}\left(
\mathrm{sc}^{([1],2)}_{s}\left(
\mathfrak{P}^{(2)}
\right)\right)\right)\right)\right)
\tag{3}
\\
&=
\mathrm{ip}^{(2,[1])\sharp}_{s}\left(
\mathrm{CH}^{[1]}_{s}\left(
\mathrm{ip}^{([1],X)@}_{s}\left(
\mathrm{sc}^{0[\mathbf{PT}_{\boldsymbol{\mathcal{A}}}]}_{s}\left(
\mathrm{sc}^{([1],2)}_{s}\left(
\mathfrak{P}^{(2)}
\right)\right)\right)\right)\right)
\tag{4}
\\
&=
\mathrm{ip}^{(2,[1])\sharp}_{s}\left(
\mathrm{sc}^{0[\mathbf{PT}_{\boldsymbol{\mathcal{A}}}]}_{s}\left(
\mathrm{sc}^{([1],2)}_{s}\left(
\mathfrak{P}^{(2)}
\right)\right)\right).
\tag{5}
\end{align*}
\end{flushleft}

In the just stated chain of equalities, the first equality unpacks the definition of the $0$-source operation just defined; the second equality unpacks the definition of the mappings $\mathrm{ip}^{(2,0)\sharp}$ and $\mathrm{sc}^{(0,2)}$, according to Definition~\ref{DDScTgZ}; the third equality follows from Proposition~\ref{PCHBasicEq} and Proposition~\ref{PPthCatAlg}; the fourth equality follows from Claim~\ref{CIsoIpfc}; finally, the last equality follows from Theorem~\ref{TIso}.

Thus, it follows that $\mathrm{sc}^{0\mathbf{Pth}_{\boldsymbol{\mathcal{A}}^{(2)}}}_{s}(\mathfrak{P}^{(2)})$ is a $(2,[1])$-identity second-order path associated to the path term class $\mathrm{sc}^{0[\mathbf{PT}_{\boldsymbol{\mathcal{A}}}]}_{s}(
\mathrm{sc}^{([1],2)}_{s}(
\mathfrak{P}^{(2)}
))$. In this regard note that, according to Corollary~\ref{CDScTgZII}, we have that 
\[
\mathrm{sc}^{0[\mathbf{PT}_{\boldsymbol{\mathcal{A}}}]}_{s}\left(
\mathrm{sc}^{([1],2)}_{s}\left(
\mathfrak{P}^{(2)}
\right)\right)
=
\mathrm{sc}^{0[\mathbf{PT}_{\boldsymbol{\mathcal{A}}}]}_{s}\left(
\mathrm{tg}^{([1],2)}_{s}\left(
\mathfrak{P}^{(2)}
\right)\right).
\]

Claim~\ref{CDPthCatAlgScZ} follows.

\textsf{(5)} For every $s\in S$, the operation of $0$-target $\mathrm{tg}_{s}^{0\mathbf{Pth}^{(1,2)}_{\boldsymbol{\mathcal{A}}^{(2)}}}$, abbreviated to $\mathrm{tg}_{s}^{0\mathbf{Pth}_{\boldsymbol{\mathcal{A}}^{(2)}}}$,  from $\mathrm{Pth}_{\boldsymbol{\mathcal{A}}^{(2)},s}$ to $\mathrm{Pth}_{\boldsymbol{\mathcal{A}}^{(2)},s}$, associated to $\mathrm{tg}^{0}_{s}$ assigns to a second-order path $\mathfrak{P}^{(2)}\in\mathrm{Pth}_{\boldsymbol{\mathcal{A}}^{(2)},s}$
the second-order $\mathbf{c}$-path  given by
\[
\mathrm{tg}_{s}^{0\mathbf{Pth}_{\boldsymbol{\mathcal{A}}^{(2)}}}\left(
\mathfrak{P}^{(2)}\right)=
\mathrm{ip}^{(2,0)\sharp}_{s}\left(
\mathrm{tg}^{(0,2)}_{s}\left(\mathfrak{P}^{(2)}
\right)\right).
\]

\begin{claim}\label{CDPthCatAlgTgZ} Let $s$ be a sort in $S$, $\mathbf{c}$ a word in $S^{\star}$, and  $\mathfrak{P}^{(2)}$ a second-order $\mathbf{c}$-path in $\mathrm{Pth}_{\boldsymbol{\mathcal{A}}^{(2)},s}$, then 
$\mathrm{tg}^{0
\mathbf{Pth}_{\boldsymbol{\mathcal{A}}^{(2)}}
}_{s}(\mathfrak{P}^{(2)})$ is the $(2,[1])$-identity second-order path in $\mathrm{Pth}_{\boldsymbol{\mathcal{A}}^{(2)},s}$ given by
\allowdisplaybreaks
\begin{align*}
\mathrm{tg}^{0
\mathbf{Pth}_{\boldsymbol{\mathcal{A}}^{(2)}}
}_{s}\left(\mathfrak{P}^{(2)}\right)
&=
\mathrm{ip}^{(2,[1])\sharp}_{s}\left(
\mathrm{tg}^{0[\mathbf{PT}_{\boldsymbol{\mathcal{A}}}]}_{s}\left(
\mathrm{sc}^{([1],2)}_{s}\left(
\mathfrak{P}^{(2)}
\right)\right)\right)
\\&=
\mathrm{ip}^{(2,[1])\sharp}_{s}\left(
\mathrm{tg}^{0[\mathbf{PT}_{\boldsymbol{\mathcal{A}}}]}_{s}\left(
\mathrm{tg}^{([1],2)}_{s}\left(
\mathfrak{P}^{(2)}
\right)\right)\right).
\end{align*}
\end{claim}

The proof is similar to that of Claim~\ref{CDPthCatAlgScZ}.

\textsf{(6)} For every $s\in S$, the partial binary operation of $0$-composition $\circ_{s}^{0\mathbf{Pth}^{(1,2)}_{\boldsymbol{\mathcal{A}}^{(2)}}}$, abbreviated to $\circ_{s}^{0\mathbf{Pth}_{\boldsymbol{\mathcal{A}}^{(2)}}}$,  from $\mathrm{Pth}_{\boldsymbol{\mathcal{A}}^{(2)},ss}$ to $\mathrm{Pth}_{\boldsymbol{\mathcal{A}}^{(2)},s}$, associated to $\circ^{0}_{s}$ assigns to a pair of second-order paths $\mathfrak{P}^{(2)}$ and $\mathfrak{Q}^{(2)}$ in $\mathrm{Pth}_{\boldsymbol{\mathcal{A}}^{(2)},s}$, where, for a unique word $\mathbf{c}\in S^{\star}$, $\mathfrak{P}^{(2)}$ is a second-order $\mathbf{c}$-path in $\boldsymbol{\mathcal{A}}^{(2)}$  of the form
$$
\mathfrak{P}^{(2)}
=
\left(
([P_{i}]_{s})_{i\in\bb{\mathbf{c}}+1},
(\mathfrak{p}^{(2)}_{i})_{i\in\bb{\mathbf{c}}},
(T^{(1)}_{i})_{i\in\bb{\mathbf{c}}}
\right),
$$
and, for a unique $\mathbf{d}\in S^{\star}$, $\mathfrak{Q}^{(2)}$ is a second-order $\mathbf{d}$-path in $\boldsymbol{\mathcal{A}}^{(2)}$  of the form
$$
\mathfrak{Q}^{(2)}
=
\left(
([Q_{j}]_{s})_{j\in\bb{\mathbf{d}}+1},
(\mathfrak{q}^{(2)}_{j})_{j\in\bb{\mathbf{d}}},
(U^{(1)}_{j})_{j\in\bb{\mathbf{d}}}
\right),
$$ 
such that 
$$\mathrm{sc}^{(0,2)}_{s}\left(\mathfrak{Q}^{(2)}\right)=\mathrm{tg}^{(0,2)}_{s}\left(\mathfrak{P}^{(2)}\right),$$ 
i.e., if we find ourselves in the situation depicted in Figure~\ref{FDPthCatAlgCompZ},  
the second-order $\mathbf{e}$-path in $\boldsymbol{\mathcal{A}}^{(2)}$ of sort $s$ is given precisely by
$$
\mathfrak{Q}^{(2)}
\circ^{0
\mathbf{Pth}_{\boldsymbol{\mathcal{A}}^{(2)}}
}_{s}\mathfrak{P}^{(2)}=
\left(
([R_{k}]_{s})_{k\in\bb{\mathbf{e}}+1},
(\mathfrak{r}^{(2)}_{k})_{k\in\bb{\mathbf{e}}}
(V^{(1)}_{k})_{k\in\bb{\mathbf{e}}}
\right),
$$

\begin{figure}
\begin{center}
\begin{tikzpicture}
[ACliment/.style={-{To [angle'=45, length=5.75pt, width=4pt, round]}}
, scale=0.75, 
AClimentD/.style={double equal sign distance,
-implies
}
]

\node[] (0) at (0,0) [] {$\mathrm{sc}^{(0,2)}_{s}(\mathfrak{P}^{(2)})$};
\node[] (1) at (5,0) [] {$\mathrm{tg}^{(0,2)}_{s}(\mathfrak{P}^{(2)})$};
\node[] (a) at (2.5,.7) [] {};
\node[] (b) at (2.5,-.7) [] {};

\node[] () at (6.5,-.1) [] {$=$};

\node[] (0p) at (8,0) [] {$\mathrm{sc}^{(0,2)}_{s}(\mathfrak{Q}^{(2)})$};
\node[] (1p) at (13,0) [] {$\mathrm{tg}^{(0,2)}_{s}(\mathfrak{Q}^{(2)})$};
\node[] (ap) at (10.5,.7) [] {};
\node[] (bp) at (10.5,-.7) [] {};

\draw[ACliment, bend left] (0) to node [above] {$\mathrm{sc}^{([1],2)}_{s}(\mathfrak{P}^{(2)})$} (1); 
\draw[ACliment, bend right] (0) to node [below] {$ \mathrm{tg}^{([1],2)}_{s}(\mathfrak{P}^{(2)})$} (1); 
\draw[AClimentD] (a) to node [right] {$ \mathfrak{P}^{(2)}$} (b); 

\draw[ACliment, bend left] (0p) to node [above] {$\mathrm{sc}^{([1],2)}_{s}(\mathfrak{Q}^{(2)})$} (1p); 
\draw[ACliment, bend right] (0p) to node [below] {$ \mathrm{tg}^{([1],2)}_{s}(\mathfrak{Q}^{(2)})$} (1p); 
\draw[AClimentD] (ap) to node [right] {$ \mathfrak{Q}^{(2)}$} (bp); 
\end{tikzpicture}
\end{center}
\caption{Requirement for the $0$-composition of two second-order paths.}\label{FDPthCatAlgCompZ}
\end{figure}

where $\mathbf{e}=\mathbf{c}\curlywedge \mathbf{d}$, and
\begin{alignat*}{2}
[R_{k}]_{s} &=
\begin{cases}
[Q_{0}
\circ^{0\mathbf{PT}_{\boldsymbol{\mathcal{A}}}}_{s}
P_{k}]_{s},
\\
[Q_{k-\bb{\mathbf{c}}}
\circ^{0\mathbf{PT}_{\boldsymbol{\mathcal{A}}}}_{s}
P_{\bb{\mathbf{c}}}]_{s},
\end{cases} 
&\,&\begin{array}{l}
\text{if $k\in \bb{\mathbf{c}}+1$;} \\
\text{if $k\in [\bb{\mathbf{c}}+1,\bb{\mathbf{e}}+1]$,}
\end{array}
\\
\mathfrak{r}^{(2)}_{k} &=
\begin{cases}
\mathfrak{p}^{(2)}_{k},\\
\mathfrak{q}^{(2)}_{k-\bb{\mathbf{c}}},\\
\end{cases}
&\,&\begin{array}{l}
\text{if $k\in \bb{\mathbf{c}}$;} \\
\text{if $k\in [\bb{\mathbf{c}},\bb{\mathbf{e}}]$,}
\end{array}
\\
V^{(1)}_{k} &=
\begin{cases}
\mathrm{CH}^{(1)}_{s}(
\mathrm{ip}^{(1,X)@}_{s}(
Q_{0}
))
\circ_{s}^{0\mathbf{PT}_{\boldsymbol{\mathcal{A}}}}T^{(1)}_{k},\\
U^{(1)}_{k-\bb{\mathbf{c}}}
\circ_{s}^{0\mathbf{PT}_{\boldsymbol{\mathcal{A}}}}
\mathrm{CH}^{(1)}_{s}(
\mathrm{ip}^{(1,X)@}_{s}(
P_{\bb{\mathbf{c}}}
))
,
\end{cases} 
&\,&\begin{array}{l}
\text{if $k\in \bb{\mathbf{c}}$;} \\
\text{if $k\in [\bb{\mathbf{c}},\bb{\mathbf{e}}]$,}
\end{array}
\end{alignat*}

That is, $\mathfrak{Q}^{(2)}
\circ^{0
\mathbf{Pth}_{\boldsymbol{\mathcal{A}}^{(2)}}
}_{s}\mathfrak{P}^{(2)}$ is a sequence of length $\bb{\mathbf{c}}+\bb{\mathbf{d}}$, the sum of the lengths of the original pair of second-order paths $\mathfrak{P}^{(2)}$ and $\mathfrak{Q}^{(2)}$. 

Moreover, for every $k\in\bb{\mathbf{c}}+1$, the $1$-constituent of $\mathfrak{Q}^{(2)}
\circ^{0
\mathbf{Pth}_{\boldsymbol{\mathcal{A}}^{(2)}}
}_{s}\mathfrak{P}^{(2)}$ at position $k$ is given by $[Q_{0}
\circ^{0\mathbf{PT}_{\boldsymbol{\mathcal{A}}}}_{s}
P_{k}]_{s}$.  Note that, for every $k\in\bb{\mathbf{c}}+1$, the following chain of equalities holds
\allowdisplaybreaks
\begin{align*}
\mathrm{sc}_{s}^{(0,1)}\left(
\mathrm{ip}^{(1,X)@}_{s}\left(
Q_{0}
\right)\right)
&=
\mathrm{sc}^{(0,[1])}_{s}\left(
\mathrm{ip}^{([1],X)@}_{s}\left(
\mathrm{sc}^{([1],2)}_{s}\left(
\mathfrak{Q}^{(2)}
\right)\right)\right)
\tag{1}
\\&=
\mathrm{sc}^{(0,2)}_{s}\left(
\mathfrak{Q}^{(2)}
\right)
\tag{2}
\\&=
\mathrm{tg}^{(0,2)}_{s}\left(
\mathfrak{P}^{(2)}
\right)
\tag{3}
\\&=
\mathrm{tg}^{(0,[1])}_{s}\left(
\mathrm{ip}^{([1],X)@}_{s}\left(
\mathrm{tg}^{([1],2)}_{s}\left(
\mathfrak{P}^{(2)}
\right)\right)\right)
\tag{4}
\\&=
\mathrm{tg}^{(0,1)}_{s}\left(
\mathrm{ip}^{(1,X)@}_{s}\left(
P_{\bb{\mathbf{c}}}
\right)\right)
\tag{5}
\\&=
\mathrm{tg}^{(0,1)}_{s}\left(
\mathrm{ip}^{(1,X)@}_{s}\left(
P_{k}
\right)\right)
.
\tag{6}
\end{align*}

The first equality recovers the definition of the mappings $\mathrm{sc}^{(0,[1])}$, $\mathrm{ip}^{([1],X)@}$ and $\mathrm{sc}^{([1],2)}$; the second equality unpacks the definition of the $(0,2)$-source; the third equality follows from the hypothesis, since $\mathrm{sc}^{(0,2)}_{s}(\mathfrak{Q}^{(2)})=\mathrm{tg}^{(0,2)}_{s}(\mathfrak{P}^{(2)})$; the fourth equality recovers the definition of the $(0,2)$-target of the second-order path $\mathfrak{P}^{(2)}$; the fifth equality unpacks the definition of the $([1],2)$-target; finally, the sixth equality follows from Proposition~\ref{PDScTgZ}.

Therefore, for every $k\in\bb{\mathbf{c}}+1$, from Proposition~\ref{PPTCatAlg}, we have that $Q_{0}
\circ^{0\mathbf{PT}_{\boldsymbol{\mathcal{A}}}}_{s}
P_{k}$ is a path term in $\mathrm{PT}_{\boldsymbol{\mathcal{A}}, s}$ and, thus,  $[Q_{0}
\circ^{0\mathbf{PT}_{\boldsymbol{\mathcal{A}}}}_{s}
P_{k}]_{s}$ is a path term class in $[\mathrm{PT}_{\boldsymbol{\mathcal{A}}}]_{s}$.

On the other hand, for every $k\in [\bb{\mathbf{c}}+1, \bb{\mathbf{e}}+1]$,  the $1$-constituent of $\mathfrak{Q}^{(2)}
\circ^{0
\mathbf{Pth}_{\boldsymbol{\mathcal{A}}^{(2)}}
}_{s}\mathfrak{P}^{(2)}$ at position $k$ is given by $[Q_{k-\bb{\mathbf{c}}}
\circ^{0\mathbf{PT}_{\boldsymbol{\mathcal{A}}}}_{s}
P_{\bb{\mathbf{c}}}]_{s}$.  Note that, for every $k\in [\bb{\mathbf{c}}+1, \bb{\mathbf{e}}+1]$, the following chain of equalities holds
\allowdisplaybreaks
\begin{align*}
\mathrm{sc}_{s}^{(0,1)}\left(
\mathrm{ip}^{(1,X)@}_{s}\left(
Q_{k-\bb{\mathbf{c}}}
\right)\right)
&=
\mathrm{sc}_{s}^{(0,1)}\left(
\mathrm{ip}^{(1,X)@}_{s}\left(
Q_{0}
\right)\right)
\tag{1}
\\&=
\mathrm{sc}^{(0,[1])}_{s}\left(
\mathrm{ip}^{([1],X)@}_{s}\left(
\mathrm{sc}^{([1],2)}_{s}\left(
\mathfrak{Q}^{(2)}
\right)\right)\right)
\tag{2}
\\&=
\mathrm{sc}^{(0,2)}_{s}\left(
\mathfrak{Q}^{(2)}
\right)
\tag{3}
\\&=
\mathrm{tg}^{(0,2)}_{s}\left(
\mathfrak{P}^{(2)}
\right)
\tag{4}
\\&=
\mathrm{tg}^{(0,[1])}_{s}\left(
\mathrm{ip}^{([1],X)@}_{s}\left(
\mathrm{tg}^{([1],2)}_{s}\left(
\mathfrak{P}^{(2)}
\right)\right)\right)
\tag{5}
\\&=
\mathrm{tg}^{(0,1)}_{s}\left(
\mathrm{ip}^{(1,X)@}_{s}\left(
P_{\bb{\mathbf{c}}}
\right)\right)
.
\tag{6}
\end{align*}

The first equality follows from Proposition~\ref{PDScTgZ}; the second equality recovers the definition of the mappings $\mathrm{sc}^{(0,[1])}$, $\mathrm{ip}^{([1],X)@}$ and $\mathrm{sc}^{([1],2)}$; the third equation recovers the $(0,2)$-target of the second-order path $\mathfrak{Q}^{(2)}$; the fourth equality follows from the hypothesis, since $\mathrm{sc}^{(0,2)}_{s}(\mathfrak{Q}^{(2)})=\mathrm{tg}^{(0,2)}_{s}(\mathfrak{P}^{(2)})$; the fifth equality recovers the definition of the $(0,2)$-target of the second-order path $\mathfrak{P}^{(2)}$; finally, the last equality unpacks the definition of the $([1],2)$-target.

Therefore, for every $k\in [\bb{\mathbf{c}}+1, \bb{\mathbf{e}}+1]$, from Proposition~\ref{PPTCatAlg}, $Q_{k-\bb{\mathbf{c}}}
\circ^{0\mathbf{PT}_{\boldsymbol{\mathcal{A}}}}_{s}
P_{\bb{\mathbf{c}}}$ is a path term in $\mathrm{PT}_{\boldsymbol{\mathcal{A}}, s}$ and, thus,  $[Q_{k-\bb{\mathbf{c}}}
\circ^{0\mathbf{PT}_{\boldsymbol{\mathcal{A}}}}_{s}
P_{\bb{\mathbf{c}}}]_{s}$ is a path term class in $[\mathrm{PT}_{\boldsymbol{\mathcal{A}}}]_{s}$.

Moreover, the $\bb{\mathbf{e}}$-family of second-order rewrite rules in $\mathfrak{Q}^{(2)}
\circ^{0
\mathbf{Pth}_{\boldsymbol{\mathcal{A}}^{(2)}}
}_{s}\mathfrak{P}^{(2)}$ is the concatenation of the respective $\bb{\mathbf{c}}$-family of second-order rewrite rules in $\mathfrak{P}^{(2)}$ and the $\bb{\mathbf{d}}$-family of second-order rewrite rules in $\mathfrak{Q}^{(2)}$.

With regard to the $\bb{\mathbf{e}}$-family of first-order translations in $\mathfrak{Q}^{(2)}
\circ^{0
\mathbf{Pth}_{\boldsymbol{\mathcal{A}}^{(2)}}
}_{s}\mathfrak{P}^{(2)}$. On one hand, the first $\bb{\mathbf{c}}$ first-order translations are given, for every $k\in\bb{\mathbf{c}}$, by the $0$-composition in $\mathrm{PT}_{\boldsymbol{\mathcal{A}}, s}$ of $\mathrm{CH}^{(1)}_{s}(\mathrm{ip}^{(1,X)@}_{s}(Q_{0}))$, the normalized path term associated to the first path term class in $\mathfrak{Q}^{(2)}$, with $T^{(1)}_{k}$, the $k$-th first-order translation in $\mathfrak{P}^{(2)}$. This is a well-defined  first-order translation according to Definition~\ref{DUTrans}.

On the other hand, the last $\bb{\mathbf{d}}$ first-order translations of $\mathfrak{Q}^{(2)}
\circ^{0
\mathbf{Pth}_{\boldsymbol{\mathcal{A}}^{(2)}}
}_{s}\mathfrak{P}^{(2)}$ are given by the $0$-composition in $\mathrm{PT}_{\boldsymbol{\mathcal{A}}, s}$  of $U^{(1)}_{k-\bb{\mathbf{c}}}$, the $(k-\bb{\mathbf{c}})$-th  first-order translation in $\mathfrak{Q}^{(2)}$, with $\mathrm{CH}^{(1)}_{s}(\mathrm{ip}^{(1,X)@}_{s}(P_{\bb{\mathbf{c}}}))$, the normalized term associated to the last path term class in $\mathfrak{P}^{(2)}$.  This is a well-defined first-order translation according to Definition~\ref{DUTrans}.

The reader is advised to consult Figure~\ref{FDPthCatAlgCompZII} for a diagrammatic representation of the $0$-composition of two second-order paths.

\begin{center}
\begin{figure}
\begin{tikzpicture}
[
ACliment/.style={{To [angle'=45, length=5.75pt, width=4pt, round]}-}
, scale=.65, 
AClimentD/.style={double equal sign distance,-implies},
AClimentR/.style={{To [angle'=45, length=5.75pt, width=4pt, round]}-}
]

\begin{scope}[xscale=1]
\draw[draw=black,pattern=crosshatch dots, pattern color=blue!20!white,   xshift=3.5cm] (0:3.5cm and 3.5cm) arc (0:180:3.5cm and 3.5cm);
\draw[draw=black, fill=white, xshift=3.5cm] (0:3.5cm and 2.6cm) arc (0:180:3.5cm and 2.6cm) ;
\draw[draw=black, pattern=dots, pattern color=blue!20!white,   xshift=3.5cm, dashed] (0:3.5cm and 1.9cm) arc (0:180:3.5cm and 1.9cm);
\draw[draw=black, fill=white,   xshift=3.5cm, dashed] (0:3.5cm and 1.1cm) arc (0:180:3.5cm and 1.1cm);
\draw[draw=black,pattern=crosshatch dots, pattern color=blue!20!white,  xshift=3.5cm] (0:3.5cm and .35cm) arc (0:180:3.5cm and .35cm);
\node[] (sc) at (3.5,3) [] {$\scriptstyle\mathrm{sc}^{([1],2)}_{s}(\mathfrak{P}^{(2)})$};

\draw[draw=black,pattern=crosshatch dots, pattern color=blue!20!white,    xshift=3.5cm] (0:3.5cm and 3.5cm) arc (0:-180:3.5cm and 3.5cm);
\node[] (tg) at (3.5,-3) [] {$\scriptstyle\mathrm{tg}^{([1],2)}_{s}(\mathfrak{P}^{(2)})$};

\draw[draw=black, pattern=crosshatch dots, pattern color=blue!20!white, xshift=10.5cm] (0:3.5cm and 3.5cm) arc (0:180:3.5cm and 3.5cm);
\node[] (scp) at (10.5,3) {$\scriptstyle\mathrm{sc}^{([1],2)}_{s}(\mathfrak{Q}^{(2)})$};

\draw[draw=black,pattern=crosshatch dots, pattern color=blue!20!white,    xshift=10.5cm] (0:3.5cm and 3.5cm) arc (0:-180:3.5cm and 3.5cm);
\draw[draw=black, fill=white,   xshift=10.5cm] (0:3.5cm and 2.6cm) arc (0:-180:3.5cm and 2.6cm);
\draw[draw=black, pattern=dots, pattern color=blue!20!white,    xshift=10.5cm, dashed] (0:3.5cm and 1.9cm) arc (0:-180:3.5cm and 1.9cm);
\draw[draw=black, fill=white,   xshift=10.5cm, dashed] (0:3.5cm and 1.1cm) arc (0:-180:3.5cm and 1.1cm);
\draw[draw=black,pattern=crosshatch dots, pattern color=blue!20!white,   xshift=10.5cm] (0:3.5cm and .35cm) arc (0:-180:3.5cm and .35cm);
\node[] (tgp) at (10.5,-3) [] {$\scriptstyle\mathrm{tg}^{([1],2)}_{s}(\mathfrak{Q}^{(2)})$};


\node[] (0) at (0,0) [] {};
\node[] (1) at (7,0) [] {$\bullet$};
\node[] (2) at (14,0) [] {};
\node[] (a) at (3.5,0) {$\scriptstyle\mathrm{tg}^{([1],2)}(\mathfrak{P}^{(2)})$} (1);
\node[] (b) at (10.5,0){$\scriptstyle\mathrm{sc}^{([1],2)}(\mathfrak{Q}^{(2)})$} (2);

\node[] () at (3.5,2.2) [] {$\Downarrow$};
\node[] () at (3.5,.7) [] {$\Downarrow$};
\node[] () at (3.5,1.45) [] {$\mathfrak{P}^{(2)}$};

\node[] () at (3.5,-2.2) [] {\rotatebox{-90}{$\,=$}};
\node[] () at (3.5,-.7) [] {\rotatebox{-90}{$\,=$}};
\node[] () at (3.5,-1.45) [] {$\vdots$};

\node[] () at (10.5,-2.2) [] {$\Downarrow$};
\node[] () at (10.5,-.7) [] {$\Downarrow$};
\node[] () at (10.5,-1.45) [] {$\mathfrak{Q}^{(2)}$};

\node[] () at (10.5,2.2) [] {\rotatebox{-90}{$\,=$}};
\node[] () at (10.5,.7) [] {\rotatebox{-90}{$\,=$}};
\node[] () at (10.5,1.45) [] {$\vdots$};

\node[] (0RU) at (.12,1) [] {};
\node[] (BLU) at (6.88,1) [] {};
\draw[-, bend right=20] (sc) to (0RU);
\draw[ACliment, bend right=20] (BLU) to (sc);

\node[] (BRU) at (7.20,.7) [] {};
\node[] (1LU) at (13.80,.7) [] {};
\draw[-, bend right=20] (scp) to (BRU);
\draw[ACliment, bend right=20] (1LU) to (scp);

\node[] (0RD) at (.20,-.7) [] {};
\node[] (BLD) at (6.80,-.7) [] {};
\draw[-, bend left=20]  (tg) to (0RD);
\draw[ACliment, bend left=20] (BLD) to (tg);

\node[] (BRD) at (7.12,-1) [] {};
\node[] (1LD) at (13.88,-1) [] {};
\draw[ACliment, bend left=20] (tgp) to  (BRD);
\draw[ACliment, bend left=20] (1LD) to (tgp);

\draw[-, shorten <=.2cm]  (a) to (0);
\draw[ACliment, shorten >=.1cm]  (1) to (a);

\draw[ACliment, shorten <=.1cm] (b) to  (1);
\draw[ACliment,  shorten >=.2cm] (2) to (b);

\draw[draw=white, fill=white] (0:.35cm) arc (0:360:.35cm);
\draw[-] ([shift=(61.5:.35cm)]0,0) arc (61.5:89:.35cm);
\draw[-] ([shift=(271.5:.35cm)]0,0) arc (271.5:386:.35cm);
\node[] () at (-.9,.05) [] {$\scriptstyle \mathrm{sc}^{(0,2)}_{s}(\mathfrak{P}^{(2)})$};

\draw[draw=white, fill=white, xshift=7cm] (0:.2cm) arc (0:360:.2cm);
\draw[-] ([shift=(-33:.2cm)]7,0) arc (-33:105:.2cm);
\draw[-] ([shift=(147:.2cm)]7,0) arc (147:285:.2cm);
\node[] () at (7,0) [] {$\bullet$};

\draw[draw=white, fill=white, xshift=14cm] (0:.35cm) arc (0:360:.35cm);
\draw[-] ([shift=(241.5:.35cm)]14,0) arc (241.5:269:.35cm);
\draw[-] ([shift=(91.5:.35cm)]14,0) arc (91.5:206:.35cm);
\node[] () at (15,.06)  [] {$\scriptstyle \mathrm{tg}^{(0,2)}_{s}(\mathfrak{Q}^{(2)})$};

\node[] () at (7,4) [] {$\scriptstyle \mathrm{sc}^{([1],2)}_{s}(\mathfrak{Q}^{(2)})\circ^{0}_{s}
\mathrm{sc}^{([1],2)}_{s}(\mathfrak{P}^{(2)})$};
\node[] () at (7,-4) [] {$\scriptstyle \mathrm{tg}^{([1],2)}_{s}(\mathfrak{Q}^{(2)})\circ^{0}_{s}
\mathrm{tg}^{([1],2)}_{s}(\mathfrak{P}^{(2)})$};

\end{scope}
\end{tikzpicture}
\caption{The $0$-composition of two second-order paths.}\label{FDPthCatAlgCompZII}
\end{figure}
\end{center}

\begin{claim}\label{CDPthCatAlgCompZ} Let $s$ be a sort in $S$, $\mathbf{c},\mathbf{d}$ words in $S^{\star}$,  $\mathfrak{P}^{(2)}$ a second-order $\mathbf{c}$-path in $\mathrm{Pth}_{\boldsymbol{\mathcal{A}}^{(2)},s}$ and $\mathfrak{Q}^{(2)}$ a second-order $\mathbf{d}$-path in $\mathrm{Pth}_{\boldsymbol{\mathcal{A}}^{(2)},s}$ such that
$$
\mathrm{sc}^{(0,2)}_{s}\left(\mathfrak{Q}^{(2)}\right)
=
\mathrm{tg}^{(0,2)}_{s}\left(
\mathfrak{P}^{(2)}
\right).
$$
Then 
$\mathfrak{Q}^{(2)}
\circ^{0
\mathbf{Pth}_{\boldsymbol{\mathcal{A}}^{(2)}}
}_{s}\mathfrak{P}^{(2)}$ is a second-order $\mathbf{c}\curlywedge\mathbf{d}$-path in $\mathrm{Pth}_{\boldsymbol{\mathcal{A}}^{(2)},s}$ of the form
\allowdisplaybreaks
\begin{multline*}
\mathfrak{Q}^{(2)}
\circ^{0
\mathbf{Pth}_{\boldsymbol{\mathcal{A}}^{(2)}}
}_{s}\mathfrak{P}^{(2)}
\colon
\mathrm{sc}^{([1],2)}_{s}\left(\mathfrak{Q}^{(2)}\right)
\circ^{0
[\mathbf{PT}_{\boldsymbol{\mathcal{A}}}
]_{\Theta^{[1]}}}_{s}
\mathrm{sc}^{([1],2)}_{s}\left(\mathfrak{P}^{(2)}\right)
{\implies}
\\
\mathrm{tg}^{([1],2)}_{s}\left(\mathfrak{Q}^{(2)}\right)
\circ^{0
[\mathbf{PT}_{\boldsymbol{\mathcal{A}}}
]_{\Theta^{[1]}}}_{s}
\mathrm{tg}^{([1],2)}_{s}\left(\mathfrak{P}^{(2)}\right).
\end{multline*}
\end{claim}

The first conditions can be easily checked, since we have that
\allowdisplaybreaks
\begin{align*}
\mathrm{sc}^{([1],2)}_{s}\left(\mathfrak{Q}^{(2)}
\circ_{s}^{0\mathbf{Pth}_{\boldsymbol{\mathcal{A}}^{(2)}}}
\mathfrak{P}^{(2)}
\right)
&=
\left[Q_{0}\circ_{s}^{0\mathbf{PT}_{\boldsymbol{\mathcal{A}}}}
P_{0}
\right]_{s}
\\&=
\left[Q_{0}\right]_{s}
\circ_{s}^{0
[\mathbf{PT}_{\boldsymbol{\mathcal{A}}}]}
\left[P_{0}\right]_{s}
\\&=
\mathrm{sc}^{([1],2)}_{s}\left(\mathfrak{Q}^{(2)}\right)
\circ_{s}^{0
[\mathbf{PT}_{\boldsymbol{\mathcal{A}}}]}
\mathrm{sc}^{([1],2)}_{s}\left(\mathfrak{P}^{(2)}\right),\, \text{and}\\
\qquad
\\
\mathrm{tg}^{([1],2)}_{s}\left(
\mathfrak{Q}^{(2)}
\circ_{s}^{0\mathbf{Pth}_{\boldsymbol{\mathcal{A}}^{(2)}}}
\mathfrak{P}^{(2)}
\right)
&=
\left[
Q_{\bb{\mathbf{d}}}
\circ_{s}^{0\mathbf{PT}_{\boldsymbol{\mathcal{A}}}}
P_{\bb{\mathbf{c}}}
\right]_{s}
\\&=
\left[Q_{\bb{\mathbf{d}}}\right]_{s}
\circ_{s}^{0
[\mathbf{PT}_{\boldsymbol{\mathcal{A}}}]}
\left[
P_{
\bb{\mathbf{c}}}
\right]_{s}
\\&=
\mathrm{tg}^{([1],2)}_{s}\left(\mathfrak{Q}^{(2)}\right)
\circ_{s}^{0
[\mathbf{PT}_{\boldsymbol{\mathcal{A}}}]}
\mathrm{tg}^{([1],2)}_{s}\left(\mathfrak{P}^{(2)}\right).
\end{align*}

The different equations above simply unpack the corresponding definitions.

Now let $k\in\bb{\mathbf{e}}$, then either (1) $k\in\bb{\mathbf{c}}$ or (2) $k\in [\bb{\mathbf{c}}, \bb{\mathbf{e}}-1]$.

If (1), i.e., if $k\in\bb{\mathbf{c}}$, then the $1$-constituent $[R_{k}]_{s}$ is given by 
$$
[R_{k}]_{s}
=
\left[Q_{0}\circ^{0\mathbf{PT}_{\boldsymbol{\mathcal{A}}}}_{s}
P_{k}\right]_{s},
$$
while the first-order translation $V^{(1)}_{k}$ is given by
$$
V^{(1)}_{k}=
\mathrm{CH}^{(1)}_{s}\left(\mathrm{ip}^{(1,X)@}_{s}\left(
Q_{0}\right)\right)
\circ_{s}^{0\mathbf{PT}_{\boldsymbol{\mathcal{A}}}}
T^{(1)}_{k}.
$$

Let us recall that the second-order rewrite rule $\mathfrak{r}^{(2)}_{k}$ is equal to the second-order rewrite rule $\mathfrak{p}^{(2)}_{k}$. If $\mathfrak{p}^{(2)}_{k}$ has the form $([M_{k}]_{c_{k}}, [N_{k}]_{c_{k}})$, and taking into account that $\mathfrak{P}^{(2)}$ is a second-order path, we have that
\begin{multicols}{2}
\begin{itemize}
\item[(i)] $T^{(1)}_{k}(
M_{k}
)\in [P_{k}]_{s}$;
\item[(ii)] $T^{(1)}_{k}(
N_{k}
)\in [P_{k+1}]_{s}$.
\end{itemize}
\end{multicols}

Moreover, in virtue of Lemma~\ref{LWCong}, we have that 
\[
\mathrm{CH}^{(1)}_{s}\left(
\mathrm{ip}^{(1,X)@}_{s}\left(
Q_{0}
\right)\right)
\in [Q_{0}]_{s}.
\]

Thus, taking into account that, according to Definition~\ref{DThetaCong}, $\Theta^{[1]}$ is a $\Sigma^{\boldsymbol{\mathcal{A}}}$-congruence on $\mathbf{PT}_{\boldsymbol{\mathcal{A}}}$, we conclude that, for every $k\in\bb{\mathbf{c}}$
\begin{multicols}{2}
\begin{itemize}
\item[(i)] $V^{(1)}_{k}(
M_{k}
)\in [R_{k}]_{s}$;
\item[(ii)] $V^{(1)}_{k}(
N_{k}
)\in [R_{k+1}]_{s}$.
\end{itemize}
\end{multicols}

If (2), i.e., if $k\in[\bb{\mathbf{c}},\bb{\mathbf{e}}-1]$, then the $1$-constituent $[R_{k}]_{s}$ is given by
$$
[R_{k}]_{s}
=
\left[Q_{k-\bb{\mathbf{c}}}\circ^{0\mathbf{PT}_{\boldsymbol{\mathcal{A}}}}_{s}
P_{\bb{\mathbf{c}}}\right]_{s},
$$
while the first-order translation $V^{(1)}_{k}$ is given by
$$
V^{(1)}_{k}=
U^{(1)}_{k-\bb{\mathbf{c}}}
\circ_{s}^{0\mathbf{PT}_{\boldsymbol{\mathcal{A}}}}
\mathrm{CH}^{(1)}_{s}\left(\mathrm{ip}^{(1,X)@}_{s}\left(
P_{0}\right)\right).
$$

Let us recall that the second-order rewrite rule $\mathfrak{r}^{(2)}_{k}$ is equal to the second-order rewrite rule $\mathfrak{q}^{(2)}_{k-\bb{\mathbf{c}}}$. If $\mathfrak{q}^{(2)}_{k-\bb{\mathbf{c}}}$ has the form $([K_{k-\bb{\mathbf{c}}}]_{d_{k-\bb{\mathbf{c}}}}, [L_{k-\bb{\mathbf{c}}}]_{d_{k-\bb{\mathbf{c}}}})$, and taking into account that $\mathfrak{Q}^{(2)}$ is a second-order path, we have that
\begin{multicols}{2}
\begin{itemize}
\item[(i)] $U^{(1)}_{k-\bb{\mathbf{c}}}(
K_{k-\bb{\mathbf{c}}}
)\in [Q_{k-\bb{\mathbf{c}}}]_{s}$;
\item[(ii)] $U^{(1)}_{k-\bb{\mathbf{c}}}(
L_{k-\bb{\mathbf{c}}}
)\in [Q_{k-\bb{\mathbf{c}}+1}]_{s}$;
\end{itemize}
\end{multicols}

Moreover, in virtue of Lemma~\ref{LWCong}, we have that 
\[
\mathrm{CH}^{(1)}_{s}\left(
\mathrm{ip}^{(1,X)@}_{s}\left(
P_{\bb{\mathbf{c}}}
\right)\right)
\in [P_{\bb{\mathbf{c}}}]_{s}.
\]

Thus, taking into account that, according to Definition~\ref{DThetaCong}, $\Theta^{[1]}$ is a $\Sigma^{\boldsymbol{\mathcal{A}}}$-congruence on $\mathbf{PT}_{\boldsymbol{\mathcal{A}}}$, we conclude that, for every $k\in[\bb{\mathbf{c}},\bb{\mathbf{e}}-1]$
\begin{multicols}{2}
\begin{itemize}
\item[(i)] $V^{(1)}_{k}(
K_{k-\bb{\mathbf{c}}}
)\in [R_{k}]_{s}$;
\item[(ii)] $V^{(1)}_{k}(
L_{k-\bb{\mathbf{c}}}
)\in [R_{k+1}]_{s}$.
\end{itemize}
\end{multicols}

This finishes the proof of Claim~\ref{CDPthCatAlgCompZ}.

This finishes the proof of Proposition~\ref{PDPthCatAlg}.
\end{proof}

We conclude this section by stating some basic properties about the just defined partial $\Sigma^{\boldsymbol{\mathcal{A}}}$-algebra of second-order paths $\mathbf{Pth}^{(1,2)}_{\boldsymbol{\mathcal{A}}^{(2)}}$.  

We next show that the $0$-source and $0$-target operations, when acting on a second-order path, always retrieve $(2,0)$-identity second-order paths.

\begin{proposition}\label{PDScTgIpZ} 
Let $s$ be a sort in $S$ and $\mathfrak{P}^{(2)}$ a second-order path in $\mathrm{Pth}_{\boldsymbol{\mathcal{A}}^{(2)}}$. Then $\mathrm{sc}^{0\mathbf{Pth}_{\boldsymbol{\mathcal{A}}^{(2)}}}_{s}(\mathfrak{P}^{(2)})$ and $\mathrm{tg}^{0\mathbf{Pth}_{\boldsymbol{\mathcal{A}}^{(2)}}}_{s}(\mathfrak{P}^{(2)})$ are $(2,0)$-identity second-order paths.
\end{proposition}
\begin{proof}
We will only prove the first statement since the second one can be handled in a similar way.

In order to prove that $\mathrm{sc}^{0\mathbf{Pth}_{\boldsymbol{\mathcal{A}}^{(2)}}}_{s}(\mathfrak{P}^{(2)})$ is a $(2,0)$-identity second-order path, we will use the characterization introduced in Proposition~\ref{PDIpZ}, where it states that this is equivalent to prove that 
\begin{enumerate}
\item $\mathrm{sc}^{0\mathbf{Pth}_{\boldsymbol{\mathcal{A}}^{(2)}}}_{s}(\mathfrak{P}^{(2)})$ is a $(2,[1])$-identity second-order path; and 
\item $\mathrm{ip}^{([1],X)@}_{s}(\mathrm{sc}^{([1],2)}_{s}(\mathrm{sc}^{0\mathbf{Pth}_{\boldsymbol{\mathcal{A}}^{(2)}}}_{s}(\mathfrak{P}^{(2)})))$ is a $([1],0)$-identity path.
\end{enumerate}

Note that, in virtue of Claim~\ref{CDPthCatAlgScZ}, the following equality holds
\[
\mathrm{sc}^{0\mathbf{Pth}_{\boldsymbol{\mathcal{A}}^{(2)}}}_{s}\left(
\mathfrak{P}^{(2)}
\right)
=
\mathrm{ip}^{(2,[1])\sharp}_{s}\left(
\mathrm{sc}^{0[\mathbf{PT}_{\boldsymbol{\mathcal{A}}}]}_{s}\left(
\mathrm{sc}^{([1],2)}_{s}\left(
\mathfrak{P}^{2}
\right)\right)\right).
\]

Therefore, $\mathrm{sc}^{0\mathbf{Pth}_{\boldsymbol{\mathcal{A}}^{(2)}}}_{s}(
\mathfrak{P}^{(2)})$ is a $(2,[1])$-identity second-order path.  To simplify the next calculations assume that, for a path term $P$ in $\mathrm{PT}_{\boldsymbol{\mathcal{A}}, s}$ it is the case that
\[
\mathrm{sc}^{([1],2)}_{s}\left(
\mathfrak{P}^{(2)}
\right)
=
[P]_{s}.
\]

Taking this into account, the following chain of equalities holds
\begin{flushleft}
$\mathrm{ip}^{([1],X)@}_{s}\left(
\mathrm{sc}^{([1],2)}_{s}\left(
\mathrm{sc}^{0\mathbf{Pth}_{\boldsymbol{\mathcal{A}}^{(2)}}}_{s}\left(
\mathfrak{P}^{(2)}
\right)\right)\right)$
\allowdisplaybreaks
\begin{align*}
\qquad
&=
\mathrm{ip}^{([1],X)@}_{s}\left(
\mathrm{sc}^{0[\mathbf{PT}_{\boldsymbol{\mathcal{A}}}]}_{s}\left(
\mathrm{sc}^{([1],2)}_{s}\left(
\mathfrak{P}^{2}
\right)\right)
\right)
\tag{1}
\\&=
\mathrm{ip}^{([1],X)@}_{s}\left(
\mathrm{sc}^{0[\mathbf{PT}_{\boldsymbol{\mathcal{A}}}]}_{s}\left(
[P]_{s}
\right)
\right)
\tag{2}
\\&=
\mathrm{ip}^{([1],X)@}_{s}\left(
\left(
\left[
\mathrm{sc}^{0\mathbf{PT}_{\boldsymbol{\mathcal{A}}}}_{s}\left(
P
\right)
\right]_{s}
\right)
\right)
\tag{3}
\\&=
\left[
\mathrm{ip}^{(1,X)@}_{s}\left(
\mathrm{sc}^{0\mathbf{PT}_{\boldsymbol{\mathcal{A}}}}_{s}\left(
P
\right)
\right)
\right]_{s}
\tag{4}
\\&=
\left[
\mathrm{sc}^{0\mathbf{Pth}_{\boldsymbol{\mathcal{A}}}}_{s}\left(
\mathrm{ip}^{(1,X)@}_{s}\left(
P
\right)
\right)
\right]_{s}
\tag{4}
\\&=
\left[
\mathrm{ip}^{(1,0)\sharp}_{s}\left(
\mathrm{sc}^{(0,1)}_{s}\left(
\mathrm{ip}^{(1,X)@}_{s}\left(
P
\right)
\right)\right)
\right]_{s}
\tag{5}
\\&=
\mathrm{ip}^{([1],0)\sharp}_{s}\left(
\mathrm{sc}^{(0,1)}_{s}\left(
\mathrm{ip}^{(1,X)@}_{s}\left(
P
\right)
\right)\right).
\tag{6}
\end{align*}
\end{flushleft}

In the just stated chain of equalities, the first equality follows from Proposition~\ref{PDUCatHom}; the second equality follows from the description of the $([1],2)$-source of $\mathfrak{P}^{(2)}$ introduced above; the third equality unpacks the description of the $0$-source operation symbol in the many-sorted partial $\Sigma^{\boldsymbol{\mathcal{A}}}$-algebra $[\mathbf{PT}_{\boldsymbol{\mathcal{A}}}]$ according to Proposition~\ref{PPTCatAlg}; the fourth equality unpacks the description of the mapping $\mathrm{ip}^{([1],X)@}$ according to Definition~\ref{DPTQIp}; the fifth equality follows from the fact that $\mathrm{ip}^{(1,X)@}$ is a $\Sigma^{\boldsymbol{\mathcal{A}}}$-homomorphism according to Definition~\ref{DIp}; the sixth equality unravels the description of the $0$-source operation symbol in the partial $\Sigma^{\boldsymbol{\mathcal{A}}}$-algebra $\mathbf{Pth}_{\boldsymbol{\mathcal{A}}}$ according to Proposition~\ref{PPthCatAlg}; finally, the last equality recovers the definition of the $\mathrm{ip}^{([1],0)\sharp}$ mapping according to Proposition~\ref{PCHBasicEq}.

This completes the proof of Proposition~\ref{PDScTgIpZ}.
\end{proof}


\begin{remark}\label{RDCompZ} The partial operation of $0$-composition in $\mathbf{Pth}^{(1,2)}_{\boldsymbol{\mathcal{A}}^{(2)}}$ is not associative in general. Moreover, the $0$-source and $0$-target operations do not necessarily act as neutral elements for the $0$-composition.
\end{remark}

An immediate consequence of Proposition~\ref{PDPthCatAlg} is that every non-constant operation symbol in the categorial signature, regardless of its arity and coarity, when applied to a family of second-order paths not entirely composed of $(2,[1])$-identity second-order paths in its domain always return coherent head-constant echelonless second-order paths.

\begin{restatable}{corollary}{CDPthWB}
\label{CDPthWB} Let $(\mathbf{s},s)$ be an element of $(S^{\star}-\{\lambda\})\times S$, $\tau$ an operation symbol in $\Sigma^{\boldsymbol{\mathcal{A}}}_{\mathbf{s},s}-\{\mathrm{sc}^{0}_{s}, \mathrm{tg}^{0}_{s}\}$, and $(\mathfrak{P}^{(2)}_{j})_{j\in\bb{\mathbf{s}}}$ a family of second-order paths in $\mathrm{Pth}_{\boldsymbol{\mathcal{A}}^{(2)},\mathbf{s}}\cap\mathrm{Dom}(\tau)$ satisfying that, for some $j\in\bb{\mathbf{s}}$, $\mathfrak{P}^{(2)}_{j}$ is second-order path of length at least one. Then the second-order path $\tau^{\mathbf{Pth}_{\boldsymbol{\mathcal{A}}^{(2)}}}((\mathfrak{P}^{(2)}_{j})_{j\in\bb{\mathbf{s}}})$ is a  coherent head-constant echelonless second-order path.
\end{restatable}

We next state that the application of the second-order path extraction algorithm, described in Lemma~\ref{LDPthExtract}, to a second-order path, which has been obtained as the result of the action of a non-constant operation symbol on a family of second-order paths not entirely composed of $(2,[1])$-identity second-order paths in its domain, recovers its original components. In this process, however, the appearance of the second-order paths is by no means chaotic, since the operations introduced in Proposition~\ref{PDPthCatAlg} follow a defined order. We recommend the reader to take a look at the analogous process on paths as depicted in the Figures~\ref{FPthExtract} and~\ref{FPthExtractProc}.

\begin{restatable}{proposition}{PDRecov}
\label{PDRecov}
Let $(\mathbf{s},s)$ be an element of $(S^{\star}-\{\lambda\})\times S$, $\tau$ an operation symbol in $\Sigma^{\boldsymbol{\mathcal{A}}}_{\mathbf{s},s}-\{\mathrm{sc}^{0}_{s},\mathrm{tg}^{0}_{s}\}$, and $(\mathfrak{P}^{(2)}_{j})_{j\in\bb{\mathbf{s}}}$ a family of second-order paths in $\mathrm{Pth}_{\boldsymbol{\mathcal{A}}^{(2)},\mathbf{s}}\cap\mathrm{Dom}(\tau)$. If, for some $j\in\bb{\mathbf{s}}$, $\mathfrak{P}^{(2)}_{j}$ is a second-order path of length at least one, then the second-order path extraction algorithm from Lemma~\ref{LDPthExtract} applied to the second-order path $\tau^{\mathbf{Pth}_{\boldsymbol{\mathcal{A}}^{(2)}}}((\mathfrak{P}^{(2)}_{j})_{j\in\bb{\mathbf{s}}})$ retrieves the original family $(\mathfrak{P}^{(2)}_{j})_{j\in\bb{\mathbf{s}}}$.
\end{restatable}

In the following corollary we state that a one-step echelonless second-order path can always be represented as the result of the action of a univocally determined operation on the family of paths extracted from it. 

\begin{restatable}{corollary}{CDUStep}
\label{CDUStep} Let $s$ be a sort in $S$ and $\mathfrak{P}^{(2)}$ be a one-step echelonless second-order  path in $\mathrm{Pth}_{\boldsymbol{\mathcal{A}}^{(2)},s}$ of type $\tau\in \Sigma^{\boldsymbol{\mathcal{A}}}_{\mathbf{s},s}$, for a unique pair $(\mathbf{s},s)\in (S^{\star}-\{\lambda\})\times S$. Let $(\mathfrak{P}^{(2)}_{j})_{j\in \bb{\mathbf{s}}}$ be the family of second-order paths we can extract from $\mathfrak{P}^{(2)}$ in virtue of Lemma~\ref{LDPthExtract}, then 
\[
\mathfrak{P}^{(2)}=\tau^{\mathbf{Pth}_{\boldsymbol{\mathcal{A}}^{(2)}}}\left(\left(
\mathfrak{P}^{(2)}_{j}\right)_{j\in\bb{\mathbf{s}}}\right).
\] 
\end{restatable}
\begin{proof}
Let us assume that, for a sort $t\in S$, $\mathfrak{P}^{(2)}$ is a second-order $t$-path in $\mathrm{Pth}_{\boldsymbol{\mathcal{A}}^{(2)},s}$ of the form
\[
\mathfrak{P}^{(2)}=
\left(
\left(
\left[
P_{i}
\right]_{s}
\right)_{i\in 2},
\mathfrak{p}^{(2)},
T^{(1)}
\right),
\]
where $T^{(1)}$ is a first-order $t$-translation of sort $s$. Following Definition~\ref{DUTrans}, for a $n\in\mathbb{N}-\{0\}$,  $T^{(1)}$ is a, a composition of elementary first-order translations $T^{(1)}_{n-1}\circ \cdots\circ  T^{(1)}_{0}$. Let us set $T^{(1)'}=T^{(1)}_{n-2}\circ \cdots \circ T^{(1)}_{0}$.

If $T^{(1)}$ is a first-order translation associated to the operation symbol $\tau\in\Sigma^{\boldsymbol{\mathcal{A}}}_{\mathbf{s},s}$, then there exists a unique index $k\in\bb{\mathbf{s}}$, a family of paths $(\mathfrak{P}_{j})_{j\in k}\in \prod_{j\in k}\mathrm{Pth}_{\boldsymbol{\mathcal{A}},s_{j}}$, a family of paths $(\mathfrak{P}_{l})_{l\in \bb{\mathbf{s}}-(k+1)}\in \prod_{l\in \bb{\mathbf{s}}-(k+1)}\mathrm{Pth}_{\boldsymbol{\mathcal{A}},s_{l}}$  such that $s_{k}=t$ and
\begin{multline*}
T^{(1)}=\tau^{\mathbf{T}_{\Sigma^{\boldsymbol{\mathcal{A}}}}(X)}\left(
\mathrm{CH}^{(1)}_{s_{0}}\left(\mathfrak{P}_{0}\right),
\cdots,
\mathrm{CH}^{(1)}_{s_{k-1}}\left(\mathfrak{P}_{k-1}\right),
\right.
\\
\left.
T^{(1)'},
\mathrm{CH}^{(1)}_{s_{k+1}}\left(\mathfrak{P}_{k+1}\right),
\cdots,
\mathrm{CH}^{(1)}_{s_{\bb{\mathbf{s}}-1}}\left(\mathfrak{P}_{\bb{\mathbf{s}}-1}\right)
\right).
\end{multline*}

Let us note that the above description is valid both for an operation symbol $\sigma\in \Sigma_{\mathbf{s},s}$ and for the $0$-composition operation symbol $\circ^{0}_{s}$.

If $\mathfrak{p}^{(2)}=([M]_{t}, [N]_{t})$, since $\mathfrak{P}^{(2)}$ is a second-order path, we have that 
\begin{multicols}{2}
\begin{itemize}
\item[(i)] $T^{(1)}(M)\in [P_{0}]_{s}$;
\item[(ii)]  $T^{(1)}(N)\in [P_{1}]_{s}$.
\end{itemize}
\end{multicols}

By Lemma~\ref{LDPthExtract}, if $(\mathfrak{P}^{(2)})_{j\in\bb{\mathbf{s}}}$ is the family of second-order paths we can extract from $\mathfrak{P}^{(2)}$, then we have that 
\[
\mathfrak{P}^{(2)}_{j}=
\begin{cases}
\mathrm{ip}^{(2,[1])\sharp}_{s_{j}}\left(
\left[
\mathrm{CH}^{(1)}_{s_{j}}\left(
\mathfrak{P}_{j}
\right)
\right]_{s_{j}}
\right)
&\qquad\mbox{if } j\neq k;
\\
\mathfrak{P}^{(2)}_{k}
&\qquad\mbox{if } j=k,
\end{cases}
\] 
where, we recall from Lemma~\ref{LDPthExtract},
$$
\xymatrix@C=70pt{
\mathfrak{P}^{(2)}_{k}\colon
[
T^{(1)'}(M)
]_{s_{k}}
\ar@{=>}[r]^-{\text{\Small{($\mathfrak{p}^{(2)}$, $
T^{(1)'}
$)}}}
&
[
T^{(1)'}(N)
]_{s_{k}}
}
.
$$

Now, if we consider the interpretation of the operation symbol $\tau$ in the many-sorted partial $\Sigma^{\boldsymbol{\mathcal{A}}}$-algebra $\mathbf{Pth}_{\boldsymbol{\mathcal{A}}^{(2)}}$ introduced in Propositions~\ref{PDPthAlg} and~\ref{PDPthCatAlg}, we have that $\tau^{\mathbf{Pth}_{\boldsymbol{\mathcal{A}}^{(2)}}}((\mathfrak{P}^{(2)}_{j})_{j\in\bb{\mathbf{s}}})$ is a second-order $t$-path in $\mathrm{Pth}_{\boldsymbol{\mathcal{A}}^{(2)},s}$ of the form
\[
\tau^{\mathbf{Pth}_{\boldsymbol{\mathcal{A}}^{(2)}}}\left(\left(\mathfrak{P}^{(2)}_{j}
\right)_{j\in\bb{\mathbf{s}}}\right)
=
\left(
\left(
\left[Q_{i}
\right]_{s}
\right)_{i\in 2}
,
\mathfrak{p}^{(2)}
,
U^{(1)}
\right),
\]

Regarding the initial $1$-constituent of $\tau^{\mathbf{Pth}_{\boldsymbol{\mathcal{A}}^{(2)}}}((\mathfrak{P}^{(2)}_{j})_{j\in\bb{\mathbf{s}}})$,  the following chain of equalities holds
\allowdisplaybreaks
\begin{align*}
[Q_{0}]_{s}&=
\left[
\tau^{\mathbf{PT}_{\boldsymbol{\mathcal{A}}}}
\left(
\mathrm{CH}^{(1)}_{s_{0}}\left(
\mathfrak{P}_{0}
\right),
\cdots,
\mathrm{CH}^{(1)}_{s_{k-1}}\left(
\mathfrak{P}_{k-1}
\right),
\right.\right.
\\&\qquad\qquad\qquad
\left.\left.
T^{(1)'}\left(
M
\right),
\mathrm{CH}^{(1)}_{s_{k+1}}\left(
\mathfrak{P}_{k+1}
\right),
\cdots,
\mathrm{CH}^{(1)}_{s_{\bb{\mathbf{s}}-1}}\left(
\mathfrak{P}_{\bb{\mathbf{s}}-1}
\right)
\right)\right]_{s}
\tag{1}
\\&=
\left[
T^{(1)}\left(
M
\right)
\right]_{s}
\tag{2}
\\&=
\left[
P_{0}
\right]_{s}.
\tag{3}
\end{align*}

In the just stated chain of equalities, the first equality unpacks the description of the initial $1$-constituent of $\tau^{\mathbf{Pth}_{\boldsymbol{\mathcal{A}}^{(2)}}}((\mathfrak{P}^{(2)}_{j})_{j\in\bb{\mathbf{s}}})$, according to Propositions~\ref{PDPthAlg} and~\ref{PDPthCatAlg}; the second equality recovers the description of the elementary first-order translation $T^{(1)}$; finally, the last equality follows from item (i) above.

Regarding the final $1$-constituent of $\tau^{\mathbf{Pth}_{\boldsymbol{\mathcal{A}}^{(2)}}}((\mathfrak{P}^{(2)}_{j})_{j\in\bb{\mathbf{s}}})$,  the following chain of equalities holds
\allowdisplaybreaks
\begin{align*}
[Q_{1}]_{s}&=
\left[
\tau^{\mathbf{PT}_{\boldsymbol{\mathcal{A}}}}
\left(
\mathrm{CH}^{(1)}_{s_{0}}\left(
\mathfrak{P}_{0}
\right),
\cdots,
\mathrm{CH}^{(1)}_{s_{k-1}}\left(
\mathfrak{P}_{k-1}
\right),
\right.\right.
\\&\qquad\qquad\qquad
\left.\left.
T^{(1)'}\left(
N
\right),
\mathrm{CH}^{(1)}_{s_{k+1}}\left(
\mathfrak{P}_{k+1}
\right),
\cdots,
\mathrm{CH}^{(1)}_{s_{\bb{\mathbf{s}}-1}}\left(
\mathfrak{P}_{\bb{\mathbf{s}}-1}
\right)
\right)\right]_{s}
\tag{1}
\\&=
\left[
T^{(1)}\left(
N
\right)
\right]_{s}
\tag{2}
\\&=
\left[
P_{1}
\right]_{s}.
\tag{3}
\end{align*}

In the just stated chain of equalities, the first equality unpacks the description of the final $1$-constituent of $\tau^{\mathbf{Pth}_{\boldsymbol{\mathcal{A}}^{(2)}}}((\mathfrak{P}^{(2)}_{j})_{j\in\bb{\mathbf{s}}})$, according to Propositions~\ref{PDPthAlg} and~\ref{PDPthCatAlg}; the second equality recovers the description of the elementary first -order translation $T^{(1)}$; finally, the last equality follows from item (ii) above.

Finally,  $U^{(1)}$, the unique first-order translation occurring in  $\tau^{\mathbf{Pth}_{\boldsymbol{\mathcal{A}}^{(2)}}}((\mathfrak{P}^{(2)}_{j})_{j\in\bb{\mathbf{s}}})$, satisfies the following chain of equalities 
\allowdisplaybreaks
\begin{align*}
U^{(1)}&=
\tau^{\mathbf{T}_{\Sigma^{\boldsymbol{\mathcal{A}}}}(X)}\left(
\mathrm{CH}^{(1)}_{s_{0}}\left(
\mathrm{ip}^{(1,X)@}_{s_{0}}\left(
\mathrm{CH}^{(1)}_{s_{0}}\left(
\mathfrak{P}_{0}
\right)
\right)
\right),
\cdots,
\right.
\\&\qquad\qquad
\mathrm{CH}^{(1)}_{s_{k-1}}\left(
\mathrm{ip}^{(1,X)@}_{s_{j-1}}\left(
\mathrm{CH}^{(1)}_{s_{k-1}}\left(
\mathfrak{P}_{k-1}
\right)
\right)\right),
T^{(1)'},
\\&\qquad\qquad\qquad
\mathrm{CH}^{(1)}_{s_{k+1}}\left(
\mathrm{ip}^{(1,X)@}_{s_{k+1}}\left(
\mathrm{CH}^{(1)}_{s_{k+1}}\left(
\mathfrak{P}_{k+1}
\right)
\right)\right),
\cdots,
\\&\qquad\qquad\qquad\qquad\qquad
\left.
\mathrm{CH}^{(1)}_{s_{\bb{\mathbf{s}}-1}}\left(
\mathrm{ip}^{(1,X)@}_{s_{\bb{\mathbf{s}}-1}}\left(
\mathrm{CH}^{(1)}_{s_{\bb{\mathbf{s}}-1}}\left(
\mathfrak{P}_{\bb{\mathbf{s}}-1}
\right)
\right)\right)
\right)
\tag{1}
\\&=
\tau^{\mathbf{T}_{\Sigma^{\boldsymbol{\mathcal{A}}}}(X)}\left(
\mathrm{CH}^{(1)}_{s_{0}}\left(\mathfrak{P}_{0}\right),
\cdots,
\mathrm{CH}^{(1)}_{s_{k-1}}\left(\mathfrak{P}_{k-1}\right),
\right.
\\&\qquad\qquad\qquad\qquad\qquad
\left.
T^{(1)'},
\mathrm{CH}^{(1)}_{s_{k+1}}\left(\mathfrak{P}_{k+1}\right),
\cdots,
\mathrm{CH}^{(1)}_{s_{\bb{\mathbf{s}}-1}}\left(\mathfrak{P}_{\bb{\mathbf{s}}-1}\right)
\right)
\tag{2}
\\&=
T^{(1)}.\tag{3}
\end{align*}

In the just stated chain of equalities, the first equality unpacks the description of the unique first-order translation occurring in $\tau^{\mathbf{Pth}_{\boldsymbol{\mathcal{A}}^{(2)}}}((\mathfrak{P}^{(2)}_{j})_{j\in\bb{\mathbf{s}}})$, according to Propositions~\ref{PDPthAlg} and~\ref{PDPthCatAlg}; the second equality follows from Proposition~\ref{PIpCH}; finally, the last equality recovers the description of the first-order translation $T$.

All in all, we conclude that 
\[
\mathfrak{P}^{(2)}=\tau^{\mathbf{Pth}_{\boldsymbol{\mathcal{A}}^{(2)}}}\left(\left(
\mathfrak{P}^{(2)}_{j}\right)_{j\in\bb{\mathbf{s}}}\right).
\]

This completes the proof.
\end{proof}

\section{
\texorpdfstring
{On the subalgebra of identity second-order paths}
{The subalgebra of second-order identity paths}
}

In this section we study how the operations defined in $\mathbf{Pth}^{(1,2)}_{\boldsymbol{\mathcal{A}}^{(2)}}$ behave on $(2,[1])$-identity second-order paths. We first note that the interpretation of every constant in $\Sigma$ is a $(2,[1])$-identity second-order path.

\begin{remark}\label{RDConsSigma}
If, for some $s\in S$, $\sigma$ is a constant operation symbol in $\Sigma_{\lambda,s}$, then its realization as a constant in $\mathbf{Pth}^{(1,2)}_{\boldsymbol{\mathcal{A}}^{(2)}}$, i.e., $\sigma^{\mathbf{Pth}_{\boldsymbol{\mathcal{A}}^{(2)}}}$, is $([\sigma^{\mathbf{PT}_{\boldsymbol{\mathcal{A}}}}]_{s},\lambda,\lambda)$, the $(2,[1])$-identity second-order path on the path term class $[\sigma^{\mathbf{PT}_{\boldsymbol{\mathcal{A}}}}]_{s}$.
\end{remark}

We next prove that every operation on the $\Sigma$-algebra of second-order paths $\mathbf{Pth}_{\boldsymbol{\mathcal{A}}^{(2)}}^{(0,2)}$ when restricted to $(2,[1])$-identity second-order paths retrieves $(2,[1])$-identity second-order paths.

\begin{restatable}{proposition}{PDUSigma}
\label{PDUSigma} Let $(\mathbf{s},s)$ be a pair in $(S^{\star}-\{\lambda\})\times S$, $\sigma$ an operation symbol in $\Sigma_{\mathbf{s},s}$ and $(\mathfrak{P}^{(2)}_{j})_{j\in\bb{\mathbf{s}}}$ a family of $(2,[1])$-identity second-order paths in $\mathrm{Pth}_{\boldsymbol{\mathcal{A}}^{(2)},\mathbf{s}}$, where, for every $j\in\bb{\mathbf{s}}$, $\mathfrak{P}^{(2)}_{j}=\mathrm{ip}^{(2,[1])\sharp}_{s_{j}}([P_{j}]_{s_{j}})$, for a suitable family of path term classes $([P_{j}]_{s_{j}})_{j\in\bb{\mathbf{s}}}$ in $[\mathrm{PT}_{\boldsymbol{\mathcal{A}}}]_{\mathbf{s}}$. Then
$$
\sigma^{\mathbf{Pth}_{\boldsymbol{\mathcal{A}}^{(2)}}}\left(
\left(\mathfrak{P}^{(2)}_{j}\right)_{j\in\bb{\mathbf{s}}}\right)
=
\mathrm{ip}^{(2,[1])\sharp}_{s}
\left(
\left[\sigma^{\mathbf{PT}_{\boldsymbol{\mathcal{A}}}}
\left((P_{j})_{j\in\bb{\mathbf{s}}}
\right)
\right]_{s}
\right).
$$
\end{restatable}
\begin{proof}
The following chain of equalities holds
\allowdisplaybreaks
\begin{align*}
\sigma^{\mathbf{Pth}_{\boldsymbol{\mathcal{A}}^{(2)}}}\left(
\left(\mathfrak{P}^{(2)}_{j}\right)_{j\in\bb{\mathbf{s}}}\right)
&=
\sigma^{\mathbf{Pth}_{\boldsymbol{\mathcal{A}}^{(2)}}}
\left(
\left(\mathrm{ip}^{(2,[1])\sharp}_{s_{j}}
\left(
\left[P_{j}
\right]_{s_{j}}
\right)
\right)_{j\in\bb{\mathbf{s}}}
\right)\tag{1}
\\&=
\mathrm{ip}^{(2,[1])\sharp}_{s}
\left(
\left[
\sigma^{\mathbf{PT}_{\boldsymbol{\mathcal{A}}}}
\left(
(P_{j})_{j\in\bb{\mathbf{s}}}
\right)
\right]_{s}
\right).\tag{2}
\end{align*}

The first equality unpacks the description of the family $(\mathfrak{P}^{(2)}_{j})_{j\in\bb{\mathbf{s}}}$ as $(2,[1])$-identity second-order paths; finally, the last equality follows from the description of the operation $\sigma$ in the many-sorted $\Sigma$-algebra $\mathbf{Pth}^{(0,2)}_{\boldsymbol{\mathcal{A}}^{(2)}}$ introduced in Proposition~\ref{PDPthAlg}, note that, according to Claim~\ref{CDPthSigma}, $\sigma^{\mathbf{Pth}_{\boldsymbol{\mathcal{A}}^{(2)}}}(
(\mathfrak{P}^{(2)}_{j})_{j\in\bb{\mathbf{s}}})$ is a second-order path of length $0$.

This completes the proof.
\end{proof}

\begin{remark}\label{RDUSigma} Taking into account Remark~\ref{RDConsSigma} and Proposition~\ref{PDUSigma}, the $S$-sorted subset of $(2,[1])$-identity second-order paths is a closed subset of $\mathbf{Pth}_{\boldsymbol{\mathcal{A}}^{(2)}}^{(0,2)}$ and we can, therefore, consider it as a $\Sigma$-subalgebra of $\mathbf{Pth}_{\boldsymbol{\mathcal{A}}^{(2)}}^{(0,2)}$.
\end{remark}

We next note that the interpretation of every constant in $\Sigma^{\boldsymbol{\mathcal{A}}}$ associated to a rewrite rule in $\boldsymbol{\mathcal{A}}^{(2)}$ is a $(2,[1])$-identity second-order path.

\begin{restatable}{remark}{RDEchCons}
\label{RDEchCons}
If, for some $s\in S$, $\mathfrak{p}$ is a constant operation symbol in $\Sigma^{\boldsymbol{\mathcal{A}}}_{\lambda,s}$ associated to a rewrite rule $\mathfrak{p}$ in $\mathcal{A}_{s}$, then its realization as a constant in $\mathbf{Pth}^{(1,2)}_{\boldsymbol{\mathcal{A}}^{(2)}}$, i.e. $\mathfrak{p}^{\mathbf{Pth}_{\boldsymbol{\mathcal{A}}^{(2)}}}$, is $([\mathfrak{p}^{\mathbf{PT}_{\boldsymbol{\mathcal{A}}}}]_{s},\lambda,\lambda)$, the $(2,[1])$-identity second-order path on the path term class $[\mathfrak{p}^{\mathbf{PT}_{\boldsymbol{\mathcal{A}}}}]_{s}$.
\end{restatable}

We next prove that every $0$-categorial operation on the many-sorted partial $\Sigma^{\boldsymbol{\mathcal{A}}}$-algebra of second-order paths $\mathbf{Pth}_{\boldsymbol{\mathcal{A}}^{(2)}}^{(1,2)}$, when restricted to $(2,[1])$-identity second-order paths retrieves $(2,[1])$-identity second-order paths. 

\begin{restatable}{proposition}{PDUIp}
\label{PDUIp} Let $s$ be a sort in $S$ and $\mathfrak{Q}^{(2)}$, $\mathfrak{P}^{(2)}$ $(2,[1])$-identity second-order paths in $\mathrm{Pth}_{\boldsymbol{\mathcal{A}}^{(2)},s}$ of the form $\mathfrak{P}^{(2)}=\mathrm{ip}^{(2,[1])\sharp}_{s}([P]_{s})$ and $\mathfrak{Q}^{(2)}=\mathrm{ip}^{(2,[1])\sharp}_{s}([Q]_{s})$, for a suitable pair of path term classes $[P]_{s}$ and $[Q]_{s}$ in $[\mathrm{PT}_{\boldsymbol{\mathcal{A}}}]_{s}$. Then
\begin{enumerate}
\item[(i)] $\mathrm{sc}^{0\mathbf{Pth}_{\boldsymbol{\mathcal{A}}^{(2)}}}_{s}(\mathfrak{P}^{(2)})
=
\mathrm{ip}^{(2,[1])\sharp}_{s}(
[
\mathrm{sc}^{0\mathbf{PT}_{\boldsymbol{\mathcal{A}}}}_{s}(
P
)
]_{s}
)
$,
\item[(ii)] $\mathrm{tg}^{0\mathbf{Pth}_{\boldsymbol{\mathcal{A}}^{(2)}}}_{s}(\mathfrak{P}^{(2)})
=
\mathrm{ip}^{(2,[1])\sharp}_{s}(
[
\mathrm{tg}^{0\mathbf{PT}_{\boldsymbol{\mathcal{A}}}}_{s}(
P
)
]_{s}
)
$ and
\item[(iii)] if $\mathrm{sc}^{(0,2)}_{s}(\mathfrak{Q}^{(2)})=\mathrm{tg}^{(0,2)}_{s}(\mathfrak{P}^{(2)})$, then
$$
\mathfrak{Q}^{(2)}
\circ_{s}^{0\mathbf{Pth}_{\boldsymbol{\mathcal{A}}^{(2)}}}
\mathfrak{P}^{(2)}
=
\mathrm{ip}^{(2,[1])\sharp}_{s}(
[
Q
\circ^{0\mathbf{PT}_{\boldsymbol{\mathcal{A}}}}_{s}
P
]_{s}
)
.
$$
\end{enumerate}
\end{restatable} 
\begin{proof}
Items (i) and (ii) are a direct consequence of Claims~\ref{CDPthCatAlgScZ} and~\ref{CDPthCatAlgTgZ}.

Therefore, we are left with the case of $0$-composition. Assume that $\mathfrak{Q}^{(2)}$ and $\mathfrak{P}^{(2)}$ satisfy that 
$\mathrm{sc}^{(0,2)}_{s}(\mathfrak{Q}^{(2)})=
\mathrm{tg}^{(0,2)}_{s}(\mathfrak{P}^{(2)})
$.  The following chain of equalities holds
\allowdisplaybreaks
\begin{align*}
\mathfrak{Q}^{(2)}
\circ^{0\mathbf{Pth}_{\boldsymbol{\mathcal{A}}^{(2)}}}_{s}
\mathfrak{P}^{(2)}&=
\mathrm{ip}^{(2,[1])\sharp}_{s}\left(
[Q]_{s}
\right)
\circ_{s}^{0\mathbf{Pth}_{\boldsymbol{\mathcal{A}}^{(2)}}}
\mathrm{ip}^{(2,[1])\sharp}_{s}\left(
[P]_{s}
\right)\tag{1}
\\&=
\mathrm{ip}^{(2,[1])\sharp}_{s}\left(
\left[
Q
\circ_{s}^{0\mathbf{PT}_{\boldsymbol{\mathcal{A}}}}
P
\right]_{s}
\right).\tag{2}
\end{align*}

The first equality describes $\mathfrak{Q}^{(2)}$ and $\mathfrak{P}^{(2)}$ as $(2,[1])$-identity second-order paths; finally, the last  equality applies the definition of the operation of  $0$-composition on the given second-order paths according to Proposition~\ref{PDPthCatAlg}. Note that the $0$-composition of two second-order paths of length $0$ is a second-order path of length $0$.

This completes the proof of Proposition~\ref{PDUCatHom}.
\end{proof}

Several important consequences follow from the previous propositions. The first one states that the $S$-sorted subset of $(2,[1])$-identity second-order paths is a closed subset of $\mathbf{Pth}_{\boldsymbol{\mathcal{A}}^{(2)}}^{(1,2)}$ and we can, therefore, consider it as a $\Sigma^{\boldsymbol{\mathcal{A}}}$-subalgebra of $\mathbf{Pth}_{\boldsymbol{\mathcal{A}}^{(2)}}^{(1,2)}$. 

\begin{restatable}{definition}{DDUIp}
\label{DDUIp}\index{identity!second-order!$\mathrm{ip}^{(2,[1])\sharp}[[\mathrm{PT}_{\boldsymbol{\mathcal{A}}}]]$} We let $\mathrm{ip}^{(2,[1])\sharp}[[\mathrm{PT}_{\boldsymbol{\mathcal{A}}}]]$ stand for the $S$-sorted subset of $(2,[1])$-identity second-order paths of $\mathrm{Pth}_{\boldsymbol{\mathcal{A}}^{(2)}}$. By Propositions~\ref{PDUSigma},\ref{PDUIp} and Remark~\ref{RDEchCons}, we conclude that $\mathrm{ip}^{(2,[1])\sharp}[[\mathrm{PT}_{\boldsymbol{\mathcal{A}}}]]$ is a closed subset of the partial $\Sigma^{\boldsymbol{\mathcal{A}}}$-algebra $\mathbf{Pth}^{(1,2)}_{\boldsymbol{\mathcal{A}}^{(2)}}$. We let $\mathrm{ip}^{(2,[1])\sharp}[[\mathbf{PT}_{\boldsymbol{\mathcal{A}}}]]$ stand for the corresponding $\Sigma^{\boldsymbol{\mathcal{A}}}$-subalgebra of the partial $\Sigma^{\boldsymbol{\mathcal{A}}}$-algebra $\mathbf{Pth}^{(1,2)}_{\boldsymbol{\mathcal{A}}^{(2)}}$.
\end{restatable}

We next state a series of propositions aiming at proving that $\mathrm{ip}^{(2,[1])\sharp}[[\mathbf{PT}_{\boldsymbol{\mathcal{A}}}]]$ is in the variety $\mathbf{PAlg}(\boldsymbol{\mathcal{E}}^{\boldsymbol{\mathcal{A}}})$, introduced in Definition~\ref{DVar}.

We start by proving that, in the partial $\Sigma^{\boldsymbol{\mathcal{A}}}$-subalgebra of $(2,[1])$-identity second-order paths, the $0$-source and $0$-target operations are right zeroes.

\begin{proposition}\label{PDUIpVarA2} Let $s$ be a sort in $S$ and let $\mathfrak{P}^{(2)}$ be a $(2,[1])$-identity second-order path in $\mathrm{Pth}_{\boldsymbol{\mathcal{A}}^{(2)},s}$, then the following equalities hold
\begin{align*}
\mathrm{sc}^{0\mathbf{Pth}_{\boldsymbol{\mathcal{A}}^{(2)}}}_{s}\left(
\mathrm{sc}^{0\mathbf{Pth}_{\boldsymbol{\mathcal{A}}^{(2)}}}_{s}\left(
\mathfrak{P}^{(2)}
\right)
\right)
&=
\mathrm{sc}^{0\mathbf{Pth}_{\boldsymbol{\mathcal{A}}^{(2)}}}_{s}\left(
\mathfrak{P}^{(2)}
\right);
\\
\mathrm{tg}^{0\mathbf{Pth}_{\boldsymbol{\mathcal{A}}^{(2)}}}_{s}\left(
\mathrm{sc}^{0\mathbf{Pth}_{\boldsymbol{\mathcal{A}}^{(2)}}}_{s}\left(
\mathfrak{P}^{(2)}
\right)
\right)
&=
\mathrm{sc}^{0\mathbf{Pth}_{\boldsymbol{\mathcal{A}}^{(2)}}}_{s}\left(
\mathfrak{P}^{(2)}
\right);
\\
\mathrm{sc}^{0\mathbf{Pth}_{\boldsymbol{\mathcal{A}}^{(2)}}}_{s}\left(
\mathrm{tg}^{0\mathbf{Pth}_{\boldsymbol{\mathcal{A}}^{(2)}}}_{s}\left(
\mathfrak{P}^{(2)}
\right)
\right)
&=
\mathrm{tg}^{0\mathbf{Pth}_{\boldsymbol{\mathcal{A}}^{(2)}}}_{s}\left(
\mathfrak{P}^{(2)}
\right);
\\
\mathrm{tg}^{0\mathbf{Pth}_{\boldsymbol{\mathcal{A}}^{(2)}}}_{s}\left(
\mathrm{tg}^{0\mathbf{Pth}_{\boldsymbol{\mathcal{A}}^{(2)}}}_{s}\left(
\mathfrak{P}^{(2)}
\right)
\right)
&=
\mathrm{tg}^{0\mathbf{Pth}_{\boldsymbol{\mathcal{A}}^{(2)}}}_{s}\left(
\mathfrak{P}^{(2)}
\right).
\end{align*}
\end{proposition}
\begin{proof}
We will only prove the first equality since the remaining ones can be handled in a similar way.

Let us assume that the $(2,[1])$-identity second-order path is given by 
$\mathfrak{P}^{(2)}=\mathrm{ip}^{(2,[1])\sharp}_{s}([P]_{s})$ for a suitable path term class $[P]_{s}$ in $[\mathrm{PT}_{\boldsymbol{\mathcal{A}}}]_{s}$.

The following chain of equalities holds
\allowdisplaybreaks
\begin{align*}
\mathrm{sc}^{0\mathbf{Pth}_{\boldsymbol{\mathcal{A}}^{(2)}}}_{s}\left(
\mathrm{sc}^{0\mathbf{Pth}_{\boldsymbol{\mathcal{A}}^{(2)}}}_{s}\left(
\mathfrak{P}^{(2)}
\right)
\right)&=
\mathrm{sc}^{0\mathbf{Pth}_{\boldsymbol{\mathcal{A}}^{(2)}}}_{s}\left(
\mathrm{sc}^{0\mathbf{Pth}_{\boldsymbol{\mathcal{A}}^{(2)}}}_{s}\left(
\mathrm{ip}^{(2,[1])\sharp}_{s}\left(
[P]_{s}
\right)
\right)
\right)
\tag{1}
\\&=
\mathrm{ip}^{(2,[1])\sharp}_{s}\left(
\left[
\mathrm{sc}^{0\mathbf{PT}_{\boldsymbol{\mathcal{A}}}}_{s}\left(
\mathrm{sc}^{0\mathbf{PT}_{\boldsymbol{\mathcal{A}}}}_{s}\left(
P
\right)
\right)
\right]_{s}
\right)\tag{2}
\\&=
\mathrm{ip}^{(2,[1])\sharp}_{s}\left(
\left[
\mathrm{sc}^{0\mathbf{PT}_{\boldsymbol{\mathcal{A}}}}_{s}\left(
P
\right)
\right]_{s}
\right)\tag{3}
\\&=
\mathrm{sc}^{0\mathbf{Pth}_{\boldsymbol{\mathcal{A}}^{(2)}}}_{s}\left(
\mathrm{ip}^{(2,[1])\sharp}_{s}\left(
\left[
P
\right]_{s}
\right)\right)\tag{4}
\\&=
\mathrm{sc}^{0\mathbf{Pth}_{\boldsymbol{\mathcal{A}}^{(2)}}}_{s}\left(
\mathfrak{P}^{(2)}
\right).\tag{5}
\end{align*}

The first equality recovers the description of $\mathfrak{P}^{(2)}$ as a $(2,[1])$-identity second-order path; the second equality follows from Proposition~\ref{PDUIp}; the third equality follows from the fact that, by Corollary~\ref{CPTFree}, $[\mathbf{PT}_{\boldsymbol{\mathcal{A}}}]$ is a many-sorted partial $\Sigma^{\boldsymbol{\mathcal{A}}}$-algebra in the QE-variety $\mathbf{PAlg}(\boldsymbol{\mathcal{E}}^{\boldsymbol{\mathcal{A}}})$; the fourth equality follows from Proposition~\ref{PDUIp}; finally, the last equality recovers the description of $\mathfrak{P}^{(2)}$ as a $(2,[1])$-identity second-order path.

This completes the proof.
\end{proof}

We next prove that, in the partial $\Sigma^{\boldsymbol{\mathcal{A}}}$-subalgebra of $(2,[1])$-identity  second-order paths, the $0$-composition of two  $(2,[1])$-identity second-order paths is defined if, and only if, the $0$-source of one second-order path equals the $0$-target of the other  second-order path.

\begin{proposition}\label{PDUIpVarA3} Let $s$ be a sort in $S$ and $\mathfrak{Q}^{(2)}$, $\mathfrak{P}^{(2)}$ two $(2,[1])$-identity second-order paths in $\mathrm{Pth}_{\boldsymbol{\mathcal{A}}^{(2)},s}$. Then the following properties are equivalent
$$
\mathfrak{Q}^{(2)}\circ^{0\mathbf{Pth}_{\boldsymbol{\mathcal{A}}^{(2)}}}_{s}
\mathfrak{P}^{(2)}\mbox{ is defined }
\quad
\Leftrightarrow
\quad
\mathrm{sc}^{0\mathbf{Pth}_{\boldsymbol{\mathcal{A}}^{(2)}}}_{s}\left(
\mathfrak{Q}^{(2)}
\right)
=
\mathrm{tg}^{0\mathbf{Pth}_{\boldsymbol{\mathcal{A}}^{(2)}}}_{s}\left(
\mathfrak{P}^{(2)}
\right).
$$
\end{proposition}
\begin{proof}
Let us assume that the $(2,[1])$-identity second-order path are given by 
\begin{multicols}{2}
\begin{itemize}
\item[] $\mathfrak{P}^{(2)}=\mathrm{ip}^{(2,[1])\sharp}_{s}([P]_{s})$;
\item[] $\mathfrak{Q}^{(2)}=\mathrm{ip}^{(2,[1])\sharp}_{s}([Q]_{s})$,
\end{itemize}
\end{multicols}
for suitable path term classes $[Q]_{s}$ and $[P]_{s}$ in $[\mathrm{PT}_{\boldsymbol{\mathcal{A}}}]_{s}$.

The following sequence of equivalences holds
\begin{flushleft}
$\mathfrak{Q}^{(2)}\circ^{0\mathbf{Pth}_{\boldsymbol{\mathcal{A}}^{(2)}}}_{s}
\mathfrak{P}^{(2)}\mbox{ is defined }$
\allowdisplaybreaks
\begin{align*}
\,&\Leftrightarrow
\mathrm{sc}^{(0,2)}_{s}\left(\mathfrak{Q}^{(2)}\right)
=
\mathrm{tg}^{(0,2)}_{s}\left(\mathfrak{P}^{(2)}\right)
\tag{1}
\\&\Leftrightarrow
\mathrm{sc}^{(0,[1])}_{s}\left(
\mathrm{ip}^{([1],X)@}_{s}\left(
\mathrm{sc}^{([1],2)}_{s}\left(
\mathrm{ip}^{(2,[1])\sharp}_{s}\left(
[Q]_{s}
\right)\right)\right)\right)
\\&\qquad\qquad\qquad\qquad\qquad
=
\mathrm{tg}^{(0,[1])}_{s}\left(
\mathrm{ip}^{([1],X)@}_{s}\left(
\mathrm{tg}^{([1],2)}_{s}\left(
\mathrm{ip}^{(2,[1])\sharp}_{s}\left(
[P]_{s}
\right)\right)\right)\right)
\tag{2}
\\&\Leftrightarrow
\mathrm{sc}^{(0,[1])}_{s}\left(
\mathrm{ip}^{([1],X)@}_{s}\left(
[Q]_{s}
\right)\right)
=
\mathrm{tg}^{(0,[1])}_{s}\left(
\mathrm{ip}^{([1],X)@}_{s}\left(
[P]_{s}
\right)\right)
\tag{3}
\\&\Leftrightarrow
\mathrm{sc}^{(0,1)}_{s}\left(
\mathrm{ip}^{(1,X)@}_{s}\left(
Q
\right)\right)
=
\mathrm{tg}^{(0,1)}_{s}\left(
\mathrm{ip}^{(1,X)@}_{s}\left(
P
\right)\right)
\tag{4}
\\&\Leftrightarrow
\mathrm{ip}^{(0,1)\sharp}_{s}\left(
\mathrm{sc}^{(0,1)}_{s}\left(
\mathrm{ip}^{(1,X)@}_{s}\left(
Q
\right)\right)\right)
=
\mathrm{ip}^{(0,1)\sharp}_{s}\left(
\mathrm{tg}^{(0,1)}_{s}\left(
\mathrm{ip}^{(1,X)@}_{s}\left(
P
\right)\right)\right)
\tag{5}
\\&\Leftrightarrow
\mathrm{sc}^{0\mathbf{Pth}_{\boldsymbol{\mathcal{A}}}}_{s}\left(
\mathrm{ip}^{(1,X)@}_{s}\left(
Q
\right)\right)
=
\mathrm{tg}^{0\mathbf{Pth}_{\boldsymbol{\mathcal{A}}}}_{s}\left(
\mathrm{ip}^{(1,X)@}_{s}\left(
P
\right)\right)
\tag{6}
\\&\Leftrightarrow
\mathrm{ip}^{(1,X)@}_{s}\left(
\mathrm{sc}^{0\mathbf{PT}_{\boldsymbol{\mathcal{A}}}}_{s}\left(
Q
\right)\right)
=
\mathrm{ip}^{(1,X)@}_{s}\left(
\mathrm{tg}^{0\mathbf{PT}_{\boldsymbol{\mathcal{A}}}}_{s}\left(
P
\right)\right)
\tag{7}
\\&\Leftrightarrow
\mathrm{CH}^{(1)}_{s}\left(
\mathrm{ip}^{(1,X)@}_{s}\left(
\mathrm{sc}^{0\mathbf{PT}_{\boldsymbol{\mathcal{A}}}}_{s}\left(
Q
\right)\right)\right)
=
\mathrm{CH}^{(1)}_{s}\left(
\mathrm{ip}^{(1,X)@}_{s}\left(
\mathrm{tg}^{0\mathbf{PT}_{\boldsymbol{\mathcal{A}}}}_{s}\left(
P
\right)\right)\right)
\tag{8}
\\&\Leftrightarrow
\left[
\mathrm{ip}^{(1,X)@}_{s}\left(
\mathrm{sc}^{0\mathbf{PT}_{\boldsymbol{\mathcal{A}}}}_{s}\left(
Q
\right)\right)\right]_{s}
=
\left[
\mathrm{ip}^{(1,X)@}_{s}\left(
\mathrm{tg}^{0\mathbf{PT}_{\boldsymbol{\mathcal{A}}}}_{s}\left(
P
\right)\right)\right]_{s}
\tag{9}
\\&\Leftrightarrow
\mathrm{ip}^{([1],X)@}_{s}\left(
\left[
\mathrm{sc}^{0\mathbf{PT}_{\boldsymbol{\mathcal{A}}}}_{s}\left(
Q
\right)\right]_{s}\right)
=
\mathrm{ip}^{([1],X)@}_{s}\left(
\left[
\mathrm{tg}^{0\mathbf{PT}_{\boldsymbol{\mathcal{A}}}}_{s}\left(
P
\right)\right]_{s}\right)
\tag{10}
\\&\Leftrightarrow
\left[
\mathrm{sc}^{0\mathbf{PT}_{\boldsymbol{\mathcal{A}}}}_{s}\left(
Q
\right)\right]_{s}
=
\left[
\mathrm{tg}^{0\mathbf{PT}_{\boldsymbol{\mathcal{A}}}}_{s}\left(
P
\right)\right]_{s}
\tag{11}
\\&\Leftrightarrow
\mathrm{ip}^{(2,[1])\sharp}_{s}\left(
\left[
\mathrm{sc}^{0\mathbf{PT}_{\boldsymbol{\mathcal{A}}}}_{s}\left(
Q
\right)\right]_{s}\right)
=
\mathrm{ip}^{(2,[1])\sharp}_{s}\left(
\left[
\mathrm{tg}^{0\mathbf{PT}_{\boldsymbol{\mathcal{A}}}}_{s}\left(
P
\right)\right]_{s}\right)
\tag{12}
\\&\Leftrightarrow
\mathrm{sc}^{0\mathbf{Pth}_{\boldsymbol{\mathcal{A}}^{(2)}}}_{s}\left(
\mathfrak{Q}^{(2)}
\right)
=
\mathrm{tg}^{0\mathbf{Pth}_{\boldsymbol{\mathcal{A}}^{(2)}}}_{s}\left(
\mathfrak{P}^{(2)}
\right).
\tag{13}
\end{align*}
\end{flushleft}

The first equivalence follows from the definition of the $0$-composition introduced in Proposition~\ref{PDPthCatAlg}; the second equivalence simply unpacks the mappings $\mathrm{sc}^{(0,2)}$ and $\mathrm{tg}^{(0,2)}$ introduced in Definition~\ref{DDScTgZ} and takes into account the nature of $\mathfrak{Q}^{(2)}$ and $\mathfrak{P}^{(2)}$ as $(2,[1])$-identity second-order paths; the  third equivalence follows from Proposition~\ref{PBasicEq}; the fourth equivalence simply unpacks the mappings $\mathrm{sc}^{(0,[1])}$ and $\mathrm{tg}^{(0,[1])}$, introduced in Definition~\ref{DCHUZ}, and the mapping $\mathrm{ip}^{([1],X)@}$, introduced in Definition~\ref{DPTQIp}; the fifth equivalence follows from Proposition~\ref{PBasicEq}; the sixth equivalence follows from the definition of the $0$-source and $0$-target introduced in Proposition~\ref{PPthCatAlg}; the seventh equivalence follows from the fact that, by Definition~\ref{DIp}, $\mathrm{ip}^{(1,X)@}$ is a $\Sigma^{\boldsymbol{\mathcal{A}}}$-homomorphism; the eighth equivalence follows from Lemma~\ref{LThetaCong}; the ninth equivalence simply unpacks the path classes; the tenth equivalence simply unpacks the mapping $\mathrm{ip}^{([1],X)@}$, introduced in Definition~\ref{DPTQIp}; the eleventh equivalence follows from the fact that, by Theorem~\ref{TIso}, $\mathrm{ip}^{([1],X)@}$ is a bijective many-sorted mapping; the twelfth equivalence follows from Proposition~\ref{PDBasicEq}; finally, the last equivalence follows from the fact that $\mathfrak{Q}^{(2)}$ and $\mathfrak{P}^{(2)}$ are $(2,[1])$-identity second-order paths and Proposition~\ref{PDUIp}.

This completes the proof.
\end{proof}

We next prove that, in the partial $\Sigma^{\boldsymbol{\mathcal{A}}}$-subalgebra of $(2,[1])$-identity second-order paths, if a $0$-composition is defined, then its $0$-source is equal to the $0$-source of the initial second-order path and its $0$-target is equal to the $0$-target of the  final second-order path.

\begin{proposition}\label{PDUIpVarA4} Let $s$ be a sort in $S$ and $\mathfrak{Q}^{(2)}$, $\mathfrak{P}^{(2)}$ $(2,[1])$-identity second-order paths in $\mathrm{Pth}_{\boldsymbol{\mathcal{A}}^{(2)},s}$ such that 
$
\mathrm{sc}^{(0,2)}_{s}(\mathfrak{Q}^{(2)})
=
\mathrm{tg}^{(0,2)}_{s}(\mathfrak{P}^{(2)}).
$ 
Then
\begin{align*}
\mathrm{sc}^{0\mathbf{Pth}_{\boldsymbol{\mathcal{A}}^{(2)}}}_{s}\left(
\mathfrak{Q}^{(2)}\circ^{0\mathbf{Pth}_{\boldsymbol{\mathcal{A}}^{(2)}}}_{s}\mathfrak{P}^{(2)}
\right)
&=
\mathrm{sc}^{0\mathbf{Pth}_{\boldsymbol{\mathcal{A}}^{(2)}}}_{s}\left(
\mathfrak{P}^{(2)}
\right);
\\
\mathrm{tg}^{0\mathbf{Pth}_{\boldsymbol{\mathcal{A}}^{(2)}}}_{s}\left(
\mathfrak{Q}^{(2)}\circ^{0\mathbf{Pth}_{\boldsymbol{\mathcal{A}}^{(2)}}}_{s}\mathfrak{P}^{(2)}
\right)
&=
\mathrm{tg}^{0\mathbf{Pth}_{\boldsymbol{\mathcal{A}}^{(2)}}}_{s}\left(
\mathfrak{Q}^{(2)}
\right).
\end{align*}
\end{proposition}
\begin{proof}
We will only prove the first equality since the second one can be handled in a similar way.

Let us assume that the $(2,[1])$-identity second-order paths are given by 
\begin{multicols}{2}
\begin{itemize}
\item[] $\mathfrak{P}^{(2)}=\mathrm{ip}^{(2,[1])\sharp}_{s}\left([P]_{s}\right)$;
\item[] $\mathfrak{Q}^{(2)}=\mathrm{ip}^{(2,[1])\sharp}_{s}\left([Q]_{s}\right)$,
\end{itemize}
\end{multicols}
for suitable path term classes $[Q]_{s}$ and $[P]_{s}$ in $[\mathrm{PT}_{\boldsymbol{\mathcal{A}}}]_{s}$.

The following sequence of equalities holds
\begin{flushleft}
$\mathrm{sc}^{0\mathbf{Pth}_{\boldsymbol{\mathcal{A}}^{(2)}}}_{s}\left(
\mathfrak{Q}^{(2)}\circ^{0\mathbf{Pth}_{\boldsymbol{\mathcal{A}}^{(2)}}}_{s}\mathfrak{P}^{(2)}
\right)$
\allowdisplaybreaks
\begin{align*}
\qquad&=
\mathrm{sc}^{0\mathbf{Pth}_{\boldsymbol{\mathcal{A}}^{(2)}}}_{s}\left(
\mathrm{ip}^{(2,[1])\sharp}_{s}\left(
[Q]_{s}
\right)
\circ^{0\mathbf{Pth}_{\boldsymbol{\mathcal{A}}^{(2)}}}_{s}
\mathrm{ip}^{(2,[1])\sharp}_{s}\left(
[P]_{s}
\right)
\right)
\tag{1}
\\&=
\mathrm{ip}^{(2,[1])\sharp}_{s}\left(
\left[\mathrm{sc}^{0\mathbf{PT}_{\boldsymbol{\mathcal{A}}}}_{s}\left(
Q\circ^{0\mathbf{PT}_{\boldsymbol{\mathcal{A}}}}_{s} P
\right)
\right]_{s}
\right)
\tag{2}
\\&=
\mathrm{ip}^{(2,[1])\sharp}_{s}\left(
\left[\mathrm{sc}^{0\mathbf{PT}_{\boldsymbol{\mathcal{A}}}}_{s}\left(
P
\right)
\right]_{s}
\right)
\tag{3}
\\&=
\mathrm{sc}^{0\mathbf{Pth}_{\boldsymbol{\mathcal{A}}^{(2)}}}_{s}\left(
\mathrm{ip}^{(2,[1])\sharp}_{s}\left(
[P]_{s}
\right)
\right)
\tag{4}
\\&=
\mathrm{sc}^{0\mathbf{Pth}_{\boldsymbol{\mathcal{A}}^{(2)}}}_{s}\left(
\mathfrak{P}^{(2)}
\right).
\tag{5}
\end{align*}
\end{flushleft}

The first equality unpacks the characterization of $\mathfrak{Q}^{(2)}$ and $\mathfrak{P}^{(2)}$ as $(2,[1])$-identity second-order paths; the second equality follows from Proposition~\ref{PDUIp}; the third equality follows from the fact that, by  Corollary~\ref{CPTFree}, $[\mathbf{PT}_{\boldsymbol{\mathcal{A}}}]$ is a many-sorted partial $\Sigma^{\boldsymbol{\mathcal{A}}}$-algebra in the QE-variety $\mathbf{PAlg}(\boldsymbol{\mathcal{E}}^{\boldsymbol{\mathcal{A}}})$; the fourth equality follows from Proposition~\ref{PDUIp}; finally, the last equality recovers the description of $\mathfrak{P}^{(2)}$ as a $(2,[1])$-identity second-order path.

This completes the proof.
\end{proof}

We next prove that in the many-sorted partial $\Sigma^{\boldsymbol{\mathcal{A}}}$-subalgebra of $(2,[1])$-identity second-order paths, the $0$-source and $0$-target operations retrieve right and left units, respectively, for  the $0$-composition.

\begin{proposition}\label{PDUIpVarA5} 
Let $s$ be a sort in $S$ and $\mathfrak{P}^{(2)}$ a $(2,[1])$-identity second-order path in $\mathrm{Pth}_{\boldsymbol{\mathcal{A}}^{(2)},s}$. Then the following equalities hold
\begin{align*}
\mathfrak{P}^{(2)}\circ^{0\mathbf{Pth}_{\boldsymbol{\mathcal{A}}^{(2)}}}_{s}
\mathrm{sc}^{0\mathbf{Pth}_{\boldsymbol{\mathcal{A}}^{(2)}}}_{s}\left(
\mathfrak{P}^{(2)}
\right)
&=
\mathfrak{P}^{(2)};
\\
\mathrm{tg}^{0\mathbf{Pth}_{\boldsymbol{\mathcal{A}}^{(2)}}}_{s}\left(
\mathfrak{P}^{(2)}
\right)
\circ^{0\mathbf{Pth}_{\boldsymbol{\mathcal{A}}^{(2)}}}_{s}
\mathfrak{P}^{(2)}
&=
\mathfrak{P}^{(2)}.
\end{align*}
\end{proposition}
\begin{proof}
We will only prove the first equality since the second one can be handled in a similar way.

Let us assume that the $(2,[1])$-identity second-order path $\mathfrak{P}^{(2)}$ is given by 
$\mathrm{ip}^{(2,[1])\sharp}_{s}([P]_{s})$, 
for a suitable path term class  $[P]_{s}$ in $[\mathrm{PT}_{\boldsymbol{\mathcal{A}}}]_{s}$.

The following sequence of equalities holds
\begin{flushleft}
$\mathfrak{P}^{(2)}\circ^{0\mathbf{Pth}_{\boldsymbol{\mathcal{A}}^{(2)}}}_{s}
\mathrm{sc}^{0\mathbf{Pth}_{\boldsymbol{\mathcal{A}}^{(2)}}}_{s}\left(
\mathfrak{P}^{(2)}
\right)$
\begin{align*}
\qquad&=
\mathrm{ip}^{(2,[1])\sharp}_{s}\left(
[P]_{s}
\right)
\circ^{0\mathbf{Pth}_{\boldsymbol{\mathcal{A}}^{(2)}}}_{s}
\mathrm{sc}^{0\mathbf{Pth}_{\boldsymbol{\mathcal{A}}^{(2)}}}_{s}\left(
\mathrm{ip}^{(2,[1])\sharp}_{s}\left(
[P]_{s}
\right)
\right)
\tag{1}
\\&=
\mathrm{ip}^{(2,[1])\sharp}_{s}\left(
\left[
P\circ^{0\mathbf{PT}_{\boldsymbol{\mathcal{A}}}}_{s}
\mathrm{sc}^{0\mathbf{PT}_{\boldsymbol{\mathcal{A}}}}_{s}(
P)
\right]_{s}
\right)
\tag{2}
\\&=
\mathrm{ip}^{(2,[1])\sharp}_{s}\left(
\left[
P
\right]_{s}
\right)
\tag{3}
\\&=
\mathfrak{P}^{(2)}.
\tag{4}
\end{align*}
\end{flushleft}

The first equality unpacks the characterization of $\mathfrak{P}^{(2)}$ as a $(2,[1])$-identity second-order path; the second equality follows from Proposition~\ref{PDUIp}; the third equality follows from the fact that, by Corollary~\ref{CPTFree}, $[\mathbf{PT}_{\boldsymbol{\mathcal{A}}}]$ is a many-sorted partial $\Sigma^{\boldsymbol{\mathcal{A}}}$-algebra in the QE-variety $\mathbf{PAlg}(\boldsymbol{\mathcal{E}}^{\boldsymbol{\mathcal{A}}})$; finally, the last equality recovers the description of $\mathfrak{P}^{(2)}$ as a $(2,[1])$-identity second-order path.

This completes the proof.
\end{proof}

We next prove that, in the many-sorted partial $\Sigma^{\boldsymbol{\mathcal{A}}}$-subalgebra of $(2,[1])$-identity second-order paths, the $0$-composition is associative.

\begin{proposition}\label{PDUIpVarA6} Let $s$ be a sort in $S$ and $\mathfrak{R}^{(2)}$, $\mathfrak{Q}^{(2)}$ and $\mathfrak{P}^{(2)}$ $(2,[1])$-identity second-order paths in $\mathrm{Pth}_{\boldsymbol{\mathcal{A}}^{(2)},s}$ such that 
\begin{align*}
\mathrm{sc}^{(0,2)}_{s}\left(
\mathfrak{R}^{(2)}
\right)
&=
\mathrm{tg}^{(0,2)}_{s}\left(
\mathfrak{Q}^{(2)}
\right),
&
\mathrm{sc}^{(0,2)}_{s}\left(
\mathfrak{Q}^{(2)}
\right)
&=
\mathrm{tg}^{(0,2)}_{s}\left(
\mathfrak{P}^{(2)}
\right).
\end{align*}
Then the following equality holds
\begin{align*}
\mathfrak{R}^{(2)}
\circ^{0\mathbf{Pth}_{\boldsymbol{\mathcal{A}}^{(2)}}}_{s}
\left(
\mathfrak{Q}^{(2)}
\circ^{0\mathbf{Pth}_{\boldsymbol{\mathcal{A}}^{(2)}}}_{s}
\mathfrak{P}^{(2)}
\right)
&=
\left(
\mathfrak{R}^{(2)}
\circ^{0\mathbf{Pth}_{\boldsymbol{\mathcal{A}}^{(2)}}}_{s}
\mathfrak{Q}^{(2)}
\right)
\circ^{0\mathbf{Pth}_{\boldsymbol{\mathcal{A}}^{(2)}}}_{s}
\mathfrak{P}^{(2)}.
\end{align*}
\end{proposition}

\begin{proof}
Let us assume that the $(2,[1])$-identity second-order paths are given by
\begin{align*}
\mathfrak{P}^{(2)}
&=
\mathrm{ip}^{(2,[1])\sharp}_{s}\left(
[P]_{s}
\right);
&
\mathfrak{Q}^{(2)}
&=
\mathrm{ip}^{(2,[1])\sharp}_{s}\left(
[Q]_{s}
\right);
&
\mathfrak{R}^{(2)}
&=
\mathrm{ip}^{(2,[1])\sharp}_{s}\left(
[R]_{s}
\right);
\end{align*}
for suitable path term classes $[P]_{s}$, $[Q]_{s}$, and $[R]_{s}$ in $[\mathrm{PT}_{\boldsymbol{\mathcal{A}}}]_{s}$.

The following chain of equalities holds
\begin{flushleft}
$\mathfrak{R}^{(2)}
\circ^{0\mathbf{Pth}_{\boldsymbol{\mathcal{A}}^{(2)}}}_{s}
\left(
\mathfrak{Q}^{(2)}
\circ^{0\mathbf{Pth}_{\boldsymbol{\mathcal{A}}^{(2)}}}_{s}
\mathfrak{P}^{(2)}
\right)$
\begin{align*}
\qquad&=
\mathrm{ip}^{(2,[1])\sharp}_{s}\left([R]_{s}\right)
\circ^{0\mathbf{Pth}_{\boldsymbol{\mathcal{A}}^{(2)}}}_{s}
\left(
\mathrm{ip}^{(2,[1])\sharp}_{s}\left(
[Q]_{s}
\right)
\circ^{0\mathbf{Pth}_{\boldsymbol{\mathcal{A}}^{(2)}}}_{s}
\mathrm{ip}^{(2,[1])\sharp}_{s}\left(
[P]_{s}
\right)
\right)
\tag{1}
\\&=
\mathrm{ip}^{(2,[1])\sharp}_{s}\left(
\left[R
\circ^{0\mathbf{PT}_{\boldsymbol{\mathcal{A}}}}_{s}
\left(Q
\circ^{0\mathbf{PT}_{\boldsymbol{\mathcal{A}}}}_{s}
P
\right)
\right]_{s}
\right)
\tag{2}
\\
&=
\mathrm{ip}^{(2,[1])\sharp}_{s}\left(
\left[\left(R
\circ^{0\mathbf{PT}_{\boldsymbol{\mathcal{A}}}}_{s}
Q\right)
\circ^{0\mathbf{PT}_{\boldsymbol{\mathcal{A}}}}_{s}
P
\right]_{s}
\right)
\tag{3}
\\&=
\left(
\mathrm{ip}^{(2,[1])\sharp}_{s}\left(
[R]_{s}
\right)
\circ^{0\mathbf{Pth}_{\boldsymbol{\mathcal{A}}^{(2)}}}_{s}
\mathrm{ip}^{(2,[1])\sharp}_{s}\left(
[Q]_{s}
\right)
\right)
\circ^{0\mathbf{Pth}_{\boldsymbol{\mathcal{A}}^{(2)}}}_{s}
\mathrm{ip}^{(2,[1])\sharp}_{s}\left(
[P]_{s}
\right)
\tag{4}
\\&=
\left(
\mathfrak{R}^{(2)}
\circ^{0\mathbf{Pth}_{\boldsymbol{\mathcal{A}}^{(2)}}}_{s}
\mathfrak{Q}^{(2)}
\right)
\circ^{0\mathbf{Pth}_{\boldsymbol{\mathcal{A}}^{(2)}}}_{s}
\mathfrak{P}^{(2)}.
\tag{5}
\end{align*}
\end{flushleft}

The first equality unpacks the characterization of $\mathfrak{R}^{(2)}$, $\mathfrak{Q}^{(2)}$ and $\mathfrak{P}^{(2)}$ as $(2,[1])$-identity second-order paths; the second equality follows from Proposition~\ref{PDUIp}; the third equality follows from the fact that, by Corollary~\ref{CPTFree}, $[\mathbf{PT}_{\boldsymbol{\mathcal{A}}}]$ is a many-sorted partial $\Sigma^{\boldsymbol{\mathcal{A}}}$-algebra in the QE-variety $\mathbf{PAlg}(\boldsymbol{\mathcal{E}}^{\boldsymbol{\mathcal{A}}})$; the fourth equality follows from Proposition~\ref{PDUIp}; finally, the last equality recovers the description of $\mathfrak{R}^{(2)}$, $\mathfrak{Q}^{(2)}$ and $\mathfrak{P}^{(2)}$ as $(2,[1])$-identity second-order paths.

This completes the proof.
\end{proof}

We next prove that, for the partial $\Sigma^{\boldsymbol{\mathcal{A}}}$-subalgebra of $(2,[1])$-identity  second-order paths, the $0$-source and $0$-target operations are $\Sigma$-homomorphisms.

\begin{proposition}\label{PDUIpVarA7} Let $(\mathbf{s},s)$ be a pair in $S^{\star}\times S$, $\sigma$ an operation symbol in $\Sigma_{\mathbf{s},s}$ and $(\mathfrak{P}^{(2)}_{j})_{j\in\bb{\mathbf{s}}}$ a family of $(2,[1])$-identity second-order paths in $\mathrm{Pth}_{\boldsymbol{\mathcal{A}}^{(2)},\mathbf{s}}$. Then the following equalities hold
\begin{align*}
\mathrm{sc}^{0\mathbf{Pth}_{\boldsymbol{\mathcal{A}}^{(2)}}}_{s}\left(
\sigma^{\mathbf{Pth}_{\boldsymbol{\mathcal{A}}^{(2)}}}\left(
\left(\mathfrak{P}^{(2)}_{j}
\right)_{j\in\bb{\mathbf{s}}}
\right)\right)
&=
\sigma^{\mathbf{Pth}_{\boldsymbol{\mathcal{A}}^{(2)}}}\left(
\left(\mathrm{sc}^{0\mathbf{Pth}_{\boldsymbol{\mathcal{A}}^{(2)}}}_{s_{j}}\left(
\mathfrak{P}^{(2)}_{j}
\right)
\right)_{j\in\bb{\mathbf{s}}}
\right);
\\
\mathrm{tg}^{0\mathbf{Pth}_{\boldsymbol{\mathcal{A}}^{(2)}}}_{s}\left(
\sigma^{\mathbf{Pth}_{\boldsymbol{\mathcal{A}}^{(2)}}}\left(
\left(\mathfrak{P}^{(2)}_{j}
\right)_{j\in\bb{\mathbf{s}}}
\right)\right)
&=
\sigma^{\mathbf{Pth}_{\boldsymbol{\mathcal{A}}^{(2)}}}\left(
\left(\mathrm{tg}^{0\mathbf{Pth}_{\boldsymbol{\mathcal{A}}^{(2)}}}_{s_{j}}\left(
\mathfrak{P}^{(2)}_{j}
\right)
\right)_{j\in\bb{\mathbf{s}}}
\right).
\end{align*}
\end{proposition}

\begin{proof}
We will only prove the first equality since the second one can be handled in a similar way.

Let us assume that, for every $j\in\bb{\mathbf{s}}$, the $(2,[1])$-identity second-order path $\mathfrak{P}^{(2)}_{j}$ is given by $
\mathfrak{P}^{(2)}_{j}=\mathrm{ip}^{(2,[1])\sharp}_{s_{j}}(
[P_{j}]_{s_{j}}
),
$ 
for a suitable path term class $[P_{j}]_{s_{j}}$ in $[\mathrm{PT}_{\boldsymbol{\mathcal{A}}}]_{s_{j}}$.

The following chain of equalities holds
\begin{flushleft}
$\mathrm{sc}^{0\mathbf{Pth}_{\boldsymbol{\mathcal{A}}^{(2)}}}_{s}\left(
\sigma^{\mathbf{Pth}_{\boldsymbol{\mathcal{A}}^{(2)}}}\left(\left(
\mathfrak{P}^{(2)}_{j}
\right)_{j\in\bb{\mathbf{s}}}
\right)\right)$
\allowdisplaybreaks
\begin{align*}
\qquad
&=
\mathrm{sc}^{0\mathbf{Pth}_{\boldsymbol{\mathcal{A}}^{(2)}}}_{s}\left(
\sigma^{\mathbf{Pth}_{\boldsymbol{\mathcal{A}}^{(2)}}}\left(\left(
\mathrm{ip}^{(2,[1])\sharp}_{s_{j}}\left(
[
P_{j}
]_{s_{j}}
\right)\right)_{j\in\bb{\mathbf{s}}}
\right)\right)
\tag{1}
\\&=
\mathrm{ip}^{(2,[1])\sharp}_{s}\left(
\left[
\mathrm{sc}^{0\mathbf{PT}_{\boldsymbol{\mathcal{A}}}}_{s}\left(
\sigma^{\mathbf{PT}_{\boldsymbol{\mathcal{A}}}}\left(
(P_{j}
)_{j\in\bb{\mathbf{s}}}
\right)\right)
\right]_{s}
\right)
\tag{2}
\\&=
\mathrm{ip}^{(2,[1])\sharp}_{s}\left(
\left[
\sigma^{\mathbf{PT}_{\boldsymbol{\mathcal{A}}}}\left(
\mathrm{sc}^{0\mathbf{PT}_{\boldsymbol{\mathcal{A}}}}_{s_{j}}\left(
(P_{j}
)\right)_{j\in\bb{\mathbf{s}}}
\right)
\right]_{s}
\right)
\tag{3}
\\&=
\sigma^{\mathbf{Pth}_{\boldsymbol{\mathcal{A}}^{(2)}}}\left(
\left(
\mathrm{sc}^{0\mathbf{Pth}_{\boldsymbol{\mathcal{A}}^{(2)}}}_{s_{j}}\left(
\mathrm{ip}^{(2,[1])\sharp}_{s_{j}}\left(
[P_{j}]_{j\in\bb{\mathbf{s}}}
\right)\right)\right)_{j\in\bb{\mathbf{s}}}
\right)
\tag{4}
\\&=
\sigma^{\mathbf{Pth}_{\boldsymbol{\mathcal{A}}^{(2)}}}\left(
\left(
\mathrm{sc}^{0\mathbf{Pth}_{\boldsymbol{\mathcal{A}}^{(2)}}}_{s_{j}}\left(
\mathfrak{P}^{(2)}_{j}
\right)\right)_{j\in\bb{\mathbf{s}}}
\right).
\tag{5}
\end{align*}
\end{flushleft}

The first equality applies the description of the family $(\mathfrak{P}^{(2)}_{j})_{j\in\bb{\mathbf{s}}}$ as $(2,[1])$-identity second-order paths; the second equality follows from Propositions~\ref{PDUSigma} and~\ref{PDUIp}; the third equality follows from the fact that, by Corollary~\ref{CPTFree}, $[\mathbf{PT}_{\boldsymbol{\mathcal{A}}}]$ is a partial $\Sigma^{\boldsymbol{\mathcal{A}}}$-algebra in the QE-variety $\mathbf{PAlg}(\boldsymbol{\mathcal{E}}^{\boldsymbol{\mathcal{A}}})$; the fourth equality follows from Propositions~\ref{PDUSigma} and~\ref{PDUIp}; finally, the last equality recovers the description of the family $(\mathfrak{P}^{(2)}_{j})_{j\in\bb{\mathbf{s}}}$ as $(2,[1])$-identity second-order paths.

This completes the proof.
\end{proof}

We next prove that, for the partial $\Sigma^{\boldsymbol{\mathcal{A}}}$-subalgebra of $(2,[1])$-identity  second-order paths, the operations from $\Sigma$ become homomorphisms with respect the $0$-composition operation.

\begin{proposition}\label{PDUIpVarA8} 
Let $(\mathbf{s},s)$ be a pair in $S^{\star}\times S$ and $(\mathfrak{Q}^{(2)}_{j})_{j\in\bb{\mathbf{s}}}$, $(\mathfrak{P}^{(2)}_{j})_{j\in\bb{\mathbf{s}}}$ two families of $(2,[1])$-identity second-order paths in $\mathrm{Pth}_{\boldsymbol{\mathcal{A}}^{(2)},\mathbf{s}}$ such that, for every $j\in\bb{\mathbf{s}}$,
$
\mathrm{sc}^{(0,2)}_{s_{j}}(
\mathfrak{Q}^{(2)}_{j}
)
=
\mathrm{tg}^{(0,2)}_{s_{j}}(
\mathfrak{P}^{(2)}_{j}
).
$ 
Then the following equality holds
\begin{multline*}
\sigma^{\mathbf{Pth}_{\boldsymbol{\mathcal{A}}^{(2)}}}\left(
\left(\mathfrak{Q}^{(2)}_{j}
\circ^{0\mathbf{Pth}_{\boldsymbol{\mathcal{A}}^{(2)}}}_{s_{j}}
\mathfrak{P}^{(2)}_{j}
\right)_{j\in\bb{\mathbf{s}}}
\right)
\\=
\sigma^{\mathbf{Pth}_{\boldsymbol{\mathcal{A}}^{(2)}}}\left(
\left(\mathfrak{Q}^{(2)}_{j}
\right)_{j\in\bb{\mathbf{s}}}
\right)
\circ^{0\mathbf{Pth}_{\boldsymbol{\mathcal{A}}^{(2)}}}_{s}
\sigma^{\mathbf{Pth}_{\boldsymbol{\mathcal{A}}^{(2)}}}\left(
\left(\mathfrak{P}^{(2)}_{j}
\right)_{j\in\bb{\mathbf{s}}}
\right).
\end{multline*}
\end{proposition}

\begin{proof}
Let us assume that, for every $j\in\bb{\mathbf{s}}$, the $(2,[1])$-identity second-order paths $\mathfrak{Q}^{(2)}_{j}$ and $\mathfrak{P}^{(2)}_{j}$ are given by
\begin{multicols}{2}
\begin{enumerate}
\item[] $\mathfrak{Q}^{(2)}_{j}=\mathrm{ip}^{(2,[1])\sharp}_{s_{j}}\left(
\left[Q_{j}\right]_{s_{j}}
\right),$
\item[] $\mathfrak{P}^{(2)}_{j}=\mathrm{ip}^{(2,[1])\sharp}_{s_{j}}\left(
\left[P_{j}\right]_{s_{j}}
\right),$
\end{enumerate}
\end{multicols}
for suitable path term classes $[Q_{j}]_{s_{j}}$ and $[P_{j}]_{s_{j}}$ in $[\mathrm{PT}_{\boldsymbol{\mathcal{A}}}]_{s_{j}}$.

The following chain of equalities holds
\begin{flushleft}
$
\sigma^{\mathbf{Pth}_{\boldsymbol{\mathcal{A}}^{(2)}}}\left(
\left(\mathfrak{Q}^{(2)}_{j}
\circ^{0\mathbf{Pth}_{\boldsymbol{\mathcal{A}}^{(2)}}}_{s_{j}}
\mathfrak{P}^{(2)}_{j}
\right)_{j\in\bb{\mathbf{s}}}
\right)
$
\allowdisplaybreaks
\begin{align*}
&=
\sigma^{\mathbf{Pth}_{\boldsymbol{\mathcal{A}}^{(2)}}}\left(
\left(
\mathrm{ip}^{(2,[1])\sharp}_{s_{j}}\left(
\left[Q_{j}
\right]_{s_{j}}
\right)
\circ^{0\mathbf{Pth}_{\boldsymbol{\mathcal{A}}^{(2)}}}_{s_{j}}
\mathrm{ip}^{(2,[1])\sharp}_{s_{j}}\left(
\left[P_{j}\right]_{s_{j}}
\right)
\right)_{j\in\bb{\mathbf{s}}}
\right)
\tag{1}
\\&=
\mathrm{ip}^{(2,[1])\sharp}_{s}\left(
\left[
\sigma^{\mathbf{PT}_{\boldsymbol{\mathcal{A}}}}\left(
\left(
Q_{j}
\circ^{0\mathbf{PT}_{\boldsymbol{\mathcal{A}}}}_{s_{j}}
P_{j}
\right)_{j\in\bb{\mathbf{s}}}
\right)
\right]_{s}
\right)
\tag{2}
\\&=
\mathrm{ip}^{(2,[1])\sharp}_{s}\left(
\left[
\sigma^{\mathbf{PT}_{\boldsymbol{\mathcal{A}}}}\left(
\left(
Q_{j}
\right)_{j\in\bb{\mathbf{s}}}
\right)
\circ^{0\mathbf{PT}_{\boldsymbol{\mathcal{A}}}}_{s}
\sigma^{\mathbf{PT}_{\boldsymbol{\mathcal{A}}}}\left(
\left(
P_{j}
\right)_{j\in\bb{\mathbf{s}}}
\right)
\right]_{s}
\right)
\tag{3}
\\&=
\sigma^{\mathbf{Pth}_{\boldsymbol{\mathcal{A}}^{(2)}}}\left(
\left(
\mathrm{ip}^{(2,[1])\sharp}_{s_{j}}\left(
\left[
Q_{j}
\right]_{s_{j}}
\right)
\right)_{j\in\bb{\mathbf{s}}}
\right)
\\&\qquad\qquad\qquad\qquad\qquad\quad\,\,\,
\circ^{0\mathbf{Pth}_{\boldsymbol{\mathcal{A}}^{(2)}}}_{s}
\sigma^{\mathbf{Pth}_{\boldsymbol{\mathcal{A}}^{(2)}}}\left(
\left(
\mathrm{ip}^{(2,[1])\sharp}_{s_{j}}\left(
\left[
P_{j}
\right]_{s_{j}}
\right)
\right)_{j\in\bb{\mathbf{s}}}
\right)
\tag{4}
\\&=
\sigma^{\mathbf{Pth}_{\boldsymbol{\mathcal{A}}^{(2)}}}\left(
\left(
\mathfrak{Q}^{(2)}_{j}
\right)_{j\in\bb{\mathbf{s}}}
\right)
\circ^{0\mathbf{Pth}_{\boldsymbol{\mathcal{A}}^{(2)}}}_{s}
\sigma^{\mathbf{Pth}_{\boldsymbol{\mathcal{A}}^{(2)}}}\left(
\left(
\mathfrak{P}^{(2)}_{j}
\right)_{j\in\bb{\mathbf{s}}}
\right).
\tag{5}
\end{align*}
\end{flushleft}

The first equality applies the description of the family $(\mathfrak{P}^{(2)}_{j})_{j\in\bb{\mathbf{s}}}$ as $(2,[1])$-identity second-order paths; the second equality follows from Propositions~\ref{PDUSigma} and~\ref{PDUIp}; the third equality follows from the fact that, by Corollary~\ref{CPTFree}, $[\mathbf{PT}_{\boldsymbol{\mathcal{A}}}]$ is a partial $\Sigma^{\boldsymbol{\mathcal{A}}}$-algebra in the QE-variety $\mathbf{PAlg}(\boldsymbol{\mathcal{E}}^{\boldsymbol{\mathcal{A}}})$; the fourth equality follows from Propositions~\ref{PDUSigma} and~\ref{PDUIp}; finally, the last equality recovers the description of the family $(\mathfrak{P}^{(2)}_{j})_{j\in\bb{\mathbf{s}}}$ as $(2,[1])$-identity second-order paths.

This completes the proof.
\end{proof}

We conclude this section by proving that the $\Sigma^{\boldsymbol{\mathcal{A}}}$-algebra of $(2,[1])$-identity second-order paths is a  partial $\Sigma^{\boldsymbol{\mathcal{A}}}$-algebra in the QE-variety $\mathbf{PAlg}(\boldsymbol{\mathcal{E}}^{\boldsymbol{\mathcal{A}}})$.

\begin{restatable}{proposition}{PDUIpAlg}
\label{PDUIpAlg} The $\Sigma^{\boldsymbol{\mathcal{A}}}$-subalgebra of $(2,[1])$-identity second-order paths, i.e. $\mathrm{ip}^{(2,[1])\sharp}[[\mathbf{PT}_{\boldsymbol{\mathcal{A}}}]]$, is a partial $\Sigma^{\boldsymbol{\mathcal{A}}}$-algebra in $\mathbf{PAlg}(\boldsymbol{\mathcal{E}}^{\boldsymbol{\mathcal{A}}})$.
\end{restatable}
\begin{proof}
We prove that all the axioms defining $\mathbf{PAlg}(\boldsymbol{\mathcal{E}}^{\boldsymbol{\mathcal{A}}})$, introduced in Definition~\ref{DVar}, are valid in 
$\mathrm{ip}^{(2,[1])\sharp}[[\mathbf{PT}_{\boldsymbol{\mathcal{A}}}]])$.

\textsf{Axiom A0.} Let $s$ be a sort in $S$ and $\mathfrak{P}^{(2)}$ a $(2,[1])$-identity second-order path in $\mathrm{Pth}_{\boldsymbol{\mathcal{A}}^{(2)},s}$. Then, according to Proposition~\ref{PDUIp},
$\mathrm{sc}^{0\mathbf{Pth}_{\boldsymbol{\mathcal{A}}^{(2)}}}_{s}(
\mathfrak{P}^{(2)}
),$ and  $\mathrm{tg}^{0\mathbf{Pth}_{\boldsymbol{\mathcal{A}}^{(2)}}}_{s}(
\mathfrak{P}^{(2)}
)$ 
are always defined and return $(2,[1])$-identity second-order paths.

\textsf{Axioms A1--A5.} These properties hold according to Propositions~\ref{PDUIpVarA2},~\ref{PDUIpVarA3},~\ref{PDUIpVarA4},~\ref{PDUIpVarA5} and~\ref{PDUIpVarA6}.

\textsf{Axiom A6.} Let $(\mathbf{s},s)$ be a pair in $S^{\star}\times S$, $\sigma$ an operation symbol in $\Sigma_{\mathbf{s},s}$ and $(\mathfrak{P}^{(2)})_{j\in\bb{\mathbf{s}}}$ a family of $(2,[1])$-identity second-order paths in $\mathrm{Pth}_{\boldsymbol{\mathcal{A}}^{(2)},\mathbf{s}}$. Then, according to Proposition~\ref{PDUSigma} and Remark~\ref{RDConsSigma},
$\sigma^{\mathbf{Pth}_{\boldsymbol{\mathcal{A}}^{(2)}}}((\mathfrak{P}^{(2)}_{j})_{j\in\bb{\mathbf{s}}})$
is always defined and returns a $(2,[1])$-identity second-order path.

\textsf{Axioms A7--A8.} These properties hold according to Propositions~\ref{PDUIpVarA7} and~\ref{PDUIpVarA8}.

\textsf{Axiom A9.} We recall from Proposition~\ref{PDPthCatAlg} that the interpretation of every rewrite rule in $\mathcal{A}$ is defined in the many-sorted partial $\Sigma^{\boldsymbol{\mathcal{A}}}$-algebra $\mathbf{Pth}^{(1,2)}_{\boldsymbol{\mathcal{A}}^{(2)}}$. In this regard, following Remark~\ref{RDEchCons}, the interpretation of every constant symbol associated to a rewrite rule in $\mathbf{Pth}_{\boldsymbol{\mathcal{A}}^{(2)}}$ is a $(2,[1])$-identity second-order path.

This completes the proof.
\end{proof}

\section{
\texorpdfstring
{Connecting layers $0$, $1$ and $2$}
{Connecting layers}
}

In this section we investigate how the different layers of (partial) algebras are interconnected.

\subsection{
\texorpdfstring
{Connecting layers $2$ and $1$}
{Connecting layers 2 and 1}
}

The structure of partial $\Sigma^{\boldsymbol{\mathcal{A}}}$-algebra that we have defined on $\mathrm{Pth}_{\boldsymbol{\mathcal{A}}^{(2)}}$ has a nice interaction with the many-sorted mappings of $([1],2)$-source $([1],2)$-target and $(2,[1])$-identity second-order path. 

We start by proving that the mappings $\mathrm{sc}^{([1],2)}$ and $\mathrm{tg}^{([1],2)}$ are $\Sigma^{\boldsymbol{\mathcal{A}}}$-homomorphisms from $\mathbf{Pth}_{\boldsymbol{\mathcal{A}}^{(2)}}^{(1,2)}$ to $[\mathbf{PT}_{\boldsymbol{\mathcal{A}}}]$.

\begin{restatable}{proposition}{PDUCatHom}
\label{PDUCatHom} The mappings $\mathrm{sc}^{([1],2)}$ and $\mathrm{tg}^{([1],2)}$ are $\Sigma^{\boldsymbol{\mathcal{A}}}$-homomorphisms from $\mathbf{Pth}^{(1,2)}_{\boldsymbol{\mathcal{A}}^{(2)}}$ to $[\mathbf{PT}_{\boldsymbol{\mathcal{A}}}]$.
\end{restatable}
\begin{proof}
Let us begin with the operations coming from the original signature $\Sigma$.
Let $(\mathbf{s},s)$ be an element of $S^{\star}\times S$, $\sigma$ an operation symbol in $\Sigma_{\mathbf{s},s}$ and $(\mathfrak{P}^{(2)}_{j})_{j\in\bb{\mathbf{s}}}$ a family of second-order paths in $\mathrm{Pth}_{\boldsymbol{\mathcal{A}}^{(2)},\mathbf{s}}$. Then, according to Claim~\ref{CDPthSigma}, we have that
\allowdisplaybreaks
\begin{align*}
\mathrm{sc}^{([1],2)}_{s}\left(
\sigma^{\mathbf{Pth}_{\boldsymbol{\mathcal{A}}^{(2)}}}
\left(\left(\mathfrak{P}^{(2)}_{j}\right)_{j\in\bb{\mathbf{s}}}\right)\right)
&=
\sigma^{[\mathbf{PT}_{\boldsymbol{\mathcal{A}}}]}
\left(\left(
\mathrm{sc}^{([1],2)}_{s_{j}}
\left(\mathfrak{P}^{(2)}_{j}\right)
\right)_{j\in\bb{\mathbf{s}}}\right),
\\
\mathrm{tg}^{([1],2)}_{s}\left(
\sigma^{\mathbf{Pth}_{\boldsymbol{\mathcal{A}}^{(2)}}}
\left(\left(\mathfrak{P}^{(2)}_{j}\right)_{j\in\bb{\mathbf{s}}}\right)\right)
&=
\sigma^{[\mathbf{PT}_{\boldsymbol{\mathcal{A}}}]}
\left(\left(
\mathrm{tg}^{([1],2)}_{s_{j}}
\left(\mathfrak{P}^{(2)}_{j}\right)
\right)_{j\in\bb{\mathbf{s}}}\right).
\end{align*}
The above statement also includes the particular case of constant symbols in $\Sigma_{\lambda,s}$ according to Remark~\ref{RDEchCons}.

For the case of constant symbols associated to rewrite rules $\mathfrak{p}$ in $\mathcal{A}_{s}$, for any sort $s\in S$, we have, according to Remark~\ref{RDEchCons}, the following equalities
$$
\mathrm{sc}^{([1],2)}_{s}\left(
\mathfrak{p}^{\mathbf{Pth}_{\boldsymbol{\mathcal{A}}^{(2)}}}
\right)=
\mathfrak{p}^{[\mathbf{PT}_{\boldsymbol{\mathcal{A}}}]}
=
\mathrm{tg}^{([1],2)}_{s}\left(
\mathfrak{p}^{\mathbf{Pth}_{\boldsymbol{\mathcal{A}}^{(2)}}}
\right).
$$

Let us consider the case of the $0$-source and $0$-target operation symbols in the categorial signature. Let $s$ be a sort in $S$ and $\mathfrak{P}^{(2)}$ a second-order path in $\mathrm{Pth}_{\boldsymbol{\mathcal{A}}^{(2)},s}$ then, according to Claims~\ref{CDPthCatAlgScZ} and~\ref{CDPthCatAlgTgZ}, we have that
\allowdisplaybreaks
\begin{align*}
\mathrm{sc}^{([1],2)}_{s}\left(
\mathrm{sc}^{0\mathbf{Pth}_{\boldsymbol{\mathcal{A}}^{(2)}}}_{s}\left(
\mathfrak{P}^{(2)}
\right)\right)
&=
\mathrm{sc}^{0
[\mathbf{PT}_{\boldsymbol{\mathcal{A}}}]
}_{s}\left(
\mathrm{sc}^{([1],2)}_{s}\left(
\mathfrak{P}^{(2)}
\right)\right);
\\
\mathrm{tg}^{([1],2)}_{s}\left(
\mathrm{sc}^{0\mathbf{Pth}_{\boldsymbol{\mathcal{A}}^{(2)}}}_{s}\left(
\mathfrak{P}^{(2)}
\right)\right)
&=
\mathrm{sc}^{0
[\mathbf{PT}_{\boldsymbol{\mathcal{A}}}]
}_{s}\left(
\mathrm{tg}^{([1],2)}_{s}\left(
\mathfrak{P}^{(2)}
\right)\right);
\\
\mathrm{sc}^{([1],2)}_{s}\left(
\mathrm{tg}^{0\mathbf{Pth}_{\boldsymbol{\mathcal{A}}^{(2)}}}_{s}\left(
\mathfrak{P}^{(2)}
\right)\right)
&=
\mathrm{tg}^{0
[\mathbf{PT}_{\boldsymbol{\mathcal{A}}}]
}_{s}\left(
\mathrm{sc}^{([1],2)}_{s}\left(
\mathfrak{P}^{(2)}
\right)\right);
\\
\mathrm{tg}^{([1],2)}_{s}\left(
\mathrm{tg}^{0\mathbf{Pth}_{\boldsymbol{\mathcal{A}}^{(2)}}}_{s}\left(
\mathfrak{P}^{(2)}
\right)\right)
&=
\mathrm{tg}^{0
[\mathbf{PT}_{\boldsymbol{\mathcal{A}}}]
}_{s}\left(
\mathrm{tg}^{([1],2)}_{s}\left(
\mathfrak{P}^{(2)}
\right)\right).
\end{align*}

Let us consider the case of the $0$-composition operation symbols in the categorial signature.
Let $s$ be a sort in $S$ and $\mathfrak{Q}^{(2)},\mathfrak{P}^{(2)}$ second-order paths in $\mathrm{Pth}_{\boldsymbol{\mathcal{A}}^{(2)},s}$ such that 
$$
\mathrm{sc}^{(0,2)}_{s}\left(\mathfrak{Q}^{(2)}\right)=
\mathrm{tg}^{(0,2)}_{s}\left(\mathfrak{P}^{(2)}\right).
$$
Then, according to Claim~\ref{CDPthCatAlgCompZ}, we have that 
\begin{align*}
\mathrm{sc}^{([1],2)}_{s}\left(
\mathfrak{Q}^{(2)}
\circ^{0\mathbf{Pth}_{\boldsymbol{\mathcal{A}}^{(2)}}}_{s}
\mathfrak{P}^{(2)}
\right)
&=
\mathrm{sc}^{([1],2)}_{s}\left(\mathfrak{Q}^{(2)}\right)
\circ^{0
[\mathbf{PT}_{\boldsymbol{\mathcal{A}}}]
}_{s}
\mathrm{sc}^{([1],2)}_{s}\left(\mathfrak{P}^{(2)}\right),
\\
\mathrm{tg}^{([1],2)}_{s}\left(
\mathfrak{Q}^{(2)}
\circ^{0\mathbf{Pth}_{\boldsymbol{\mathcal{A}}^{(2)}}}_{s}
\mathfrak{P}^{(2)}
\right)
&=
\mathrm{tg}^{([1],2)}_{s}\left(\mathfrak{Q}^{(2)}\right)
\circ^{0
[\mathbf{PT}_{\boldsymbol{\mathcal{A}}}]
}_{s}
\mathrm{tg}^{([1],2)}_{s}\left(\mathfrak{P}^{(2)}\right).
\end{align*}

This proves Proposition~\ref{PDUCatHom}.
\end{proof}

Also $\mathrm{ip}^{(2,[1])\sharp}$, the $(2,[1])$-identity path mapping, is a $\Sigma^{\boldsymbol{\mathcal{A}}}$-homomorphism from $[\mathbf{PT}_{\boldsymbol{\mathcal{A}}}]$ to $\mathbf{Pth}^{(1,2)}_{\boldsymbol{\mathcal{A}}^{(2)}}$.  In this regard, the following remark is introduced to simplify notation.

\begin{remark}\label{RDUIp} By Definition~\ref{DDUIp}, $\mathrm{ip}^{(2,[1])\sharp}[
[
\mathbf{PT}_{\boldsymbol{\mathcal{A}}}
]
]$ is a partial $\Sigma^{\boldsymbol{\mathcal{A}}}$-subalgebra of $\mathbf{Pth}_{\boldsymbol{\mathcal{A}}^{(2)}}^{(1,2)}$. Therefore, the inclusion mapping $\mathrm{in}^{\mathrm{ip}^{(2,[1])\sharp}[
[
\mathbf{PT}_{\boldsymbol{\mathcal{A}}}
]
]}$ from $\mathrm{ip}^{(2,[1])\sharp}[
[
\mathbf{PT}_{\boldsymbol{\mathcal{A}}}
]
]$ to $\mathbf{Pth}_{\boldsymbol{\mathcal{A}}^{(2)}}^{(1,2)}$ becomes a $\Sigma^{\boldsymbol{\mathcal{A}}}$-homomorphism.  To simplify the presentation, we will write $\mathrm{ip}^{(2,[1])\sharp}$ instead of $\mathrm{in}^{\mathrm{ip}^{(2,[1])\sharp}[
[
\mathbf{PT}_{\boldsymbol{\mathcal{A}}}
]
]}\circ\mathrm{ip}^{(2,[1])\sharp}$.
\end{remark}

To obtain the aforementioned result we will make use of the universal property of the $\Sigma^{\boldsymbol{\mathcal{A}}}$-algebra $[\mathbf{PT}_{\boldsymbol{\mathcal{A}}}]$ and the fact that, by Proposition~\ref{PDUIpAlg}, $\mathrm{ip}^{(2,[1])\sharp}([\mathbf{PT}_{\boldsymbol{\mathcal{A}}}])$ is a partial $\Sigma^{\boldsymbol{\mathcal{A}}}$-algebra
in the QE-variety $\mathbf{PAlg}(\boldsymbol{\mathcal{E}}^{\boldsymbol{\mathcal{A}}})$. 

\begin{restatable}{proposition}{PDUIpCatHom}
\label{PDUIpCatHom} The mapping $\mathrm{ip}^{(2,[1])\sharp}$ is a $\Sigma^{\boldsymbol{\mathcal{A}}}$-homomorphism from $[\mathbf{PT}_{\boldsymbol{\mathcal{A}}}]$ to $\mathbf{Pth}^{(1,2)}_{\boldsymbol{\mathcal{A}}^{(2)}}$.
\end{restatable}
\begin{proof}
Let us recall, from Definition~\ref{DDPth}, that $\mathrm{ip}^{(2,X)}$ is the $S$-sorted mapping from $X$ to $\mathbf{Pth}^{(1,2)}_{\boldsymbol{\mathcal{A}}^{(2)}}$ that, for every sort $s\in S$, sends $x\in X_{s}$ to $([x]_{s},\lambda,\lambda)$, the $(2,[1])$-identity second-order path on $[x]_{s}$. This $S$-sorted mapping corestricts to the subset of $(2,[1])$-identity second-order paths, i.e., $\mathrm{ip}^{(2,[1])\sharp}[[\mathrm{PT}_{\boldsymbol{\mathcal{A}}}]]$. Now, by Proposition~\ref{PDUIpAlg}, we have that $\mathrm{ip}^{(2,[1])\sharp}[[\mathbf{PT}_{\boldsymbol{\mathcal{A}}}]]$ is a partial $\Sigma^{\boldsymbol{\mathcal{A}}}$-algebra
in the QE-variety $\mathbf{PAlg}(\boldsymbol{\mathcal{E}}^{\boldsymbol{\mathcal{A}}})$. Moreover, by Corollary~\ref{CPTFree}, there exists a unique $\Sigma^{\boldsymbol{\mathcal{A}}}$-homomorphism $\mathrm{ip}^{(2,[1])\sharp}$ from $[\mathbf{PT}_{\boldsymbol{\mathcal{A}}}]$ to $\mathrm{ip}^{(2,[1])\sharp}[[\mathbf{PT}_{\boldsymbol{\mathcal{A}}}]]$ (to $\mathbf{Pth}^{(1,2)}_{\boldsymbol{\mathcal{A}}^{(2)}}$) such that
$$
\mathrm{ip}^{(2,[1])\sharp}\circ\eta^{([1],X)}=\mathrm{ip}^{(2,X)}, 
$$
i.e., such that the diagram in Figure~\ref{FDUCatHom} commutes.

Let us note that the $\Sigma^{\boldsymbol{\mathcal{A}}}$-homomorphism $\mathrm{ip}^{(2,[1])\sharp}$ sends, for every sort $s\in S$, a path term class $[P]_{s}\in[\mathrm{PT}_{\boldsymbol{\mathcal{A}}}]_{s}$ to $([P]_{s},\lambda,\lambda)$, the $(2,[1])$-identity second-order path on $[P]_{s}$.
\end{proof}

\begin{figure}
\begin{center}
\begin{tikzpicture}
[ACliment/.style={-{To [angle'=45, length=5.75pt, width=4pt, round]}},scale=1]
\node[] (xoq) at (0,0) [] {$X$};
\node[] (txoq) at (6,0) [] {$[\mathbf{PT}_{\boldsymbol{\mathcal{A}}}]$};
\node[] (p) at (6,-3) [] {$
\mathbf{Pth}^{(1,2)}_{\boldsymbol{\mathcal{A}}^{(2)}}
$};
\draw[ACliment]  (xoq) to node [above]
{$\eta^{([1],X)}$} (txoq);
\draw[ACliment, bend right=10]  (xoq) to node [below left] {$\mathrm{ip}^{(2,X)}$} (p);

\node[] (B0) at (6,-1.5)  [] {};
\draw[ACliment]  ($(B0)+(0,1.2)$) to node [above, fill=white] {
$\textstyle \mathrm{ip}^{(2,[1])\sharp}$
} ($(B0)+(0,-1.2)$);
\draw[ACliment, bend right]  ($(B0)+(.4,-1.2)$) to node [ below, fill=white] {
$\textstyle \mathrm{tg}^{([1],2)}$
} ($(B0)+(.4,1.2)$);
\draw[ACliment, bend left]  ($(B0)+(-.4,-1.2)$) to node [below, fill=white] {
$\textstyle \mathrm{sc}^{([1],2)}$
} ($(B0)+(-.4,1.2)$);

\end{tikzpicture}
\end{center}
\caption{$\Sigma^{\boldsymbol{\mathcal{A}}}$-homomorphisms relative to $X$ at layers 1 \& 2.}
\label{FDUCatHom}
\end{figure}

From the previous results, it follows that the partial $\Sigma^{\boldsymbol{\mathcal{A}}}$-algebra of path term classes is isomorphic to the many-sorted partial $\Sigma^{\boldsymbol{\mathcal{A}}}$-subalgebra of $(2,[1])$-identity second-order paths.

\begin{corollary}\label{CDUCatIso} The mapping $\mathrm{ip}^{(2,[1])\sharp}$ is a $\Sigma^{\boldsymbol{\mathcal{A}}}$-isomorphism from $[\mathbf{PT}_{\boldsymbol{\mathcal{A}}}]$ to $\mathrm{ip}^{(2,[1])\sharp}[[\mathbf{PT}_{\boldsymbol{\mathcal{A}}}]]$.
\end{corollary}

\begin{proof}
According to Proposition~\ref{PDUIpCatHom}, $\mathrm{ip}^{(2,[1])\sharp}$ is a many-sorted $\Sigma^{\boldsymbol{\mathcal{A}}}$-homomorphism from $[\mathbf{PT}_{\boldsymbol{\mathcal{A}}}]$ to $\mathbf{Pth}_{\boldsymbol{\mathcal{A}}^{(2)}}^{(1,2)}$ that corestricts to $\mathrm{ip}^{(2,[1])\sharp}[[\mathbf{PT}_{\boldsymbol{\mathcal{A}}}]]$. That $\mathrm{ip}^{(2,[1])\sharp}$ is surjective follows from the description of the codomain and that it is injective follows from Proposition~\ref{PDBasicEq}, since $\mathrm{sc}^{([1],2)}\circ\mathrm{ip}^{(2,[1])\sharp}=\mathrm{id}^{[\mathrm{PT}_{\boldsymbol{\mathcal{A}}}]}$ (analogously with the $([1],2)$-target).
\end{proof}

\begin{proposition}\label{PDUCatSection}
The following equations between $\Sigma^{\boldsymbol{\mathcal{A}}}$-homomorphisms hold
\begin{itemize}
\item[(i)] $\mathrm{sc}^{([1],2)}\circ \mathrm{ip}^{(2,[1])\sharp}
=\mathrm{id}^{[\mathbf{PT}_{\boldsymbol{\mathcal{A}}}]}$,
\item[(ii)] $\mathrm{tg}^{([1],2)}\circ \mathrm{ip}^{(2,[1])\sharp}
=\mathrm{id}^{[\mathbf{PT}_{\boldsymbol{\mathcal{A}}}]}$.
\end{itemize}
\end{proposition}
\begin{proof}
It follows from Proposition~\ref{PDBasicEq}.
\end{proof}

Thus, $[\mathbf{PT}_{\boldsymbol{\mathcal{A}}}]$ is a retract of $\mathbf{Pth}^{(1,2)}_{\boldsymbol{\mathcal{A}}^{(2)}}$, the
$\Sigma^{\boldsymbol{\mathcal{A}}}$-homomorphisms $\mathrm{sc}^{([1],2)}$ and $\mathrm{tg}^{([1],2)}$ are retractions of $\mathbf{Pth}^{(1,2)}_{\boldsymbol{\mathcal{A}}^{(2)}}$ onto $[\mathbf{PT}_{\boldsymbol{\mathcal{A}}}]$, and the $\Sigma^{\boldsymbol{\mathcal{A}}}$-homomorphism $\mathrm{ip}^{(2,[1])\sharp}$ is a section of both $\mathrm{sc}^{([1],2)}$ and $\mathrm{tg}^{([1],2)}$.

\begin{corollary}\label{CDUScTg}
Let $s$ be a sort in $S$ and $\mathfrak{P}^{(2)}$ a second-order path in $\mathrm{Pth}_{\boldsymbol{\mathcal{A}}^{(2)},s}$, then we have that
\begin{align*}
\left(\mathfrak{P}^{(2)},\mathrm{ip}^{(2,[1])\sharp}_{s}\left(\mathrm{sc}^{([1],2)}_{s}\left(\mathfrak{P}^{(2)}\right)\right)\right)
&\in\mathrm{Ker}\left(\mathrm{sc}^{([1],2)}\right)_{s},\\
\left(\mathfrak{P}^{(2)},\mathrm{ip}^{(2,[1])\sharp}_{s}\left(\mathrm{tg}^{([1],2)}_{s}\left(\mathfrak{P}^{(2)}\right)\right)\right)
&\in\mathrm{Ker}\left(\mathrm{tg}^{([1],2)}\right)_{s}.
\end{align*}
\end{corollary}

\begin{remark}
From Proposition~\ref{PDUCatHom} it follows that the kernel of $\mathrm{sc}^{([1],2)}$ is a $\Sigma^{\boldsymbol{\mathcal{A}}}$-congruence on $\mathbf{Pth}^{(1,2)}_{\boldsymbol{\mathcal{A}}^{(2)}}$. Therefore, we can consider the standard (epi, mono)-factoriza\-tion of $\mathrm{sc}^{([1],2)}$ given by
\begin{enumerate}
\item the projection $\mathrm{pr}^{\mathrm{Ker}(\mathrm{sc}^{([1],2)})}$ from $\mathbf{Pth}^{(1,2)}_{\boldsymbol{\mathcal{A}}^{(2)}}$ to $\mathbf{Pth}^{(1,2)}_{\boldsymbol{\mathcal{A}}^{(2)}}/{\mathrm{Ker}(\mathrm{sc}^{([1],2)})}$ that, for every $s\in S$, assigns to a second-order path $\mathfrak{P}^{(2)}$ in $\mathrm{Pth}_{\boldsymbol{\mathcal{A}}^{(2)},s}$ its equivalence class $[\mathfrak{P}^{(2)}]_{\mathrm{Ker}(\mathrm{sc}^{([1],2)})_{s}}$, and 
\item the embedding $\mathrm{sc}^{([1],2)\mathrm{m}}$ of $\mathbf{Pth}^{(1,2)}_{\boldsymbol{\mathcal{A}}^{(2)}}/{\mathrm{Ker}(\mathrm{sc}^{([1],2)})}$ into $[\mathbf{PT}_{\boldsymbol{\mathcal{A}}}]$ that, for every $s\in S$, assigns to an equivalence class $[\mathfrak{P}^{(2)}]_{\mathrm{Ker}(\mathrm{sc}^{([1],2)})_{s}}$ in the quotient $\mathrm{Pth}_{\boldsymbol{\mathcal{A}}^{(2)},s}/{\mathrm{Ker}(\mathrm{sc}^{([1],2)})_{s}}$, with $\mathfrak{P}^{(2)}\in\mathrm{Pth}_{\boldsymbol{\mathcal{A}}^{(2)},s}$, the path term class $\mathrm{sc}^{([1],2)}_{s}(\mathfrak{P}^{(2)})$, i.e., the value of the mapping $\mathrm{sc}^{([1],2)}_{s}$ at any equivalence class representative. We will refer to the mapping $\mathrm{sc}^{([1],2)\mathrm{m}}$ as the \emph{monomorphic factorization of $\mathrm{sc}^{([1],2)}$}.
\end{enumerate}

Let us recall that $(\mathrm{pr}^{\mathrm{Ker}(\mathrm{sc}^{([1],2)})}, \mathrm{sc}^{(1,2)\mathrm{m}})$ is the unique pair of $\Sigma^{\boldsymbol{\mathcal{A}}}$-homomorphisms such that $\mathrm{sc}^{([1],2)\mathrm{m}}\circ\mathrm{pr}^{\mathrm{Ker}(\mathrm{sc^{(1,2)}})}=\mathrm{sc}^{([1],2)}$, i.e., such that the upper half of the diagram in Figure~\ref{FDUIsos} commutes.

On the other hand, the kernel of $\mathrm{tg}^{([1],2)}$ is also a $\Sigma^{\boldsymbol{\mathcal{A}}}$-congruence on $\mathbf{Pth}^{(1,2)}_{\boldsymbol{\mathcal{A}}^{(2)}}$. We can consider the standard (epi, mono)-factorization of $\mathrm{tg}^{([1],2)}$ given by 
\begin{enumerate}
\item the projection $\mathrm{pr}^{\mathrm{Ker}(\mathrm{tg}^{([1],2)})}$ from $\mathbf{Pth}^{(1,2)}_{\boldsymbol{\mathcal{A}}^{(2)}}$ to $\mathbf{Pth}^{(1,2)}_{\boldsymbol{\mathcal{A}}^{(2)}}/{\mathrm{Ker}(\mathrm{tg}^{([1],2)})}$ that, for every $s\in S$, assigns to a second-order path $\mathfrak{P}^{(2)}$ in $\mathrm{Pth}_{\boldsymbol{\mathcal{A}}^{(2)},s}$ its equivalence class $[\mathfrak{P}^{(2)}]_{\mathrm{Ker}(\mathrm{tg}^{([1],2)})_{s}}$, and
\item the embedding $\mathrm{tg}^{([1],2)\mathrm{m}}$ from $\mathbf{Pth}^{(1,2)}_{\boldsymbol{\mathcal{A}}^{(2)}}/{\mathrm{Ker}(\mathrm{tg}^{([1],2)})}$ to $[\mathbf{PT}_{\boldsymbol{\mathcal{A}}}]$ that, for every $s\in S$, assigns to an equivalence class $[\mathfrak{P}^{(2)}]_{\mathrm{Ker}(\mathrm{tg}^{([1],2)})_{s}}$ in the quotient $\mathrm{Pth}_{\boldsymbol{\mathcal{A}}^{(2)},s}/{\mathrm{Ker}(\mathrm{tg}^{([1],2)})_{s}}$ and $\mathfrak{P}^{(2)}\in\mathrm{Pth}_{\boldsymbol{\mathcal{A}}^{(2)},s}$ the path term class $\mathrm{tg}^{([1],2)}_{s}(\mathfrak{P}^{(2)})$, i.e., the value of the mapping $\mathrm{tg}^{([1],2)}_{s}$ at any equivalence class representative. We will refer to the mapping $\mathrm{tg}^{(1,2)\mathrm{m}}$ as the \emph{monomorphic factorization of $\mathrm{tg}^{([1],2)}$}.
\end{enumerate}

Let us recall that $(\mathrm{pr}^{\mathrm{Ker}(\mathrm{tg}^{([1],2)})}, \mathrm{tg}^{([1],2)\mathrm{m}})$ is the unique pair of $\Sigma^{\boldsymbol{\mathcal{A}}}$-homomorphisms such that $\mathrm{tg}^{([1],2)\mathrm{m}}\circ\mathrm{pr}^{\mathrm{Ker}(\mathrm{tg}^{([1],2)})}=\mathrm{tg}^{([1],2)}$, i.e., such that the lower half of the diagram in Figure~\ref{FDUIsos} commutes.

\end{remark}

\begin{figure}
\begin{center}
\begin{tikzpicture}
[ACliment/.style={-{To [angle'=45, length=5.75pt, width=4pt, round]}}]
\node[] (txoq) at (-1,0) [] {$[\mathbf{PT}_{\boldsymbol{\mathcal{A}}}]$};
\node[] (p) at (4,3.5) [] {$\mathbf{Pth}^{(1,2)}_{\boldsymbol{\mathcal{A}}^{(2)}}$};
\node[] (p2) at (4,-3.5) [] {$\mathbf{Pth}^{(1,2)}_{\boldsymbol{\mathcal{A}}^{(2)}}$};
\node[] (pkernel) at (4,1.5) [] {$\mathbf{Pth}^{(1,2)}_{\boldsymbol{\mathcal{A}}^{(2)}}/{\mathrm{Ker}(\mathrm{sc}^{([1],2)})}$};
\node[] (pkernel2) at (4,-1.5) [] {$\mathbf{Pth}^{(1,2)}_{\boldsymbol{\mathcal{A}}^{(2)}}/{\mathrm{Ker}(\mathrm{tg}^{([1],2)})}$};
\node[] (txoq2) at (9,0) [] {$[\mathbf{PT}_{\boldsymbol{\mathcal{A}}}]$};
\draw[ACliment, bend left=15]  (txoq) to node [above left]
{$\mathrm{ip}^{(2,[1])\sharp}$} (p);
\draw[ACliment]  (txoq) to node [sloped, above] {}
(pkernel);
\draw[ACliment]  (p) to node [midway, fill=white]
{$\mathrm{pr}^{\mathrm{Ker}(\mathrm{sc}^{([1],2)})}$} (pkernel);
\draw[ACliment]  (pkernel) to node [sloped, above]
{$\mathrm{sc}^{([1],2)\mathrm{m}}$} (txoq2);
\draw[ACliment, bend left=15]  (p) to node [above right]
{$\mathrm{sc}^{([1],2)}$} (txoq2);
\draw[ACliment]  (txoq) to node  [midway, fill=white]
{$\mathrm{id}^{[\mathbf{PT}_{\boldsymbol{\mathcal{A}}}]}$} (txoq2);
\draw[ACliment, bend right=15]  (txoq) to node [below left]
{$\mathrm{ip}^{(2,[1])\sharp}$} (p2);
\draw[ACliment]  (txoq) to node [sloped, below] {}
 (pkernel2);
\draw[ACliment]  (p2) to node  [midway, fill=white]
{$\mathrm{pr}^{\mathrm{Ker}(\mathrm{tg}^{([1],2)})}$} (pkernel2);
\draw[ACliment]  (pkernel2) to node [sloped, below]
{$\mathrm{tg}^{([1],2)\mathrm{m}}$} (txoq2);
\draw[ACliment, bend right=15]  (p2) to node [below right]
{$\mathrm{tg}^{([1],2)}$} (txoq2);
\end{tikzpicture}
\end{center}
\caption{$\Sigma^{\boldsymbol{\mathcal{A}}}$-algebras and $\Sigma^{\boldsymbol{\mathcal{A}}}$-homomorphisms at layers $2$ \& $1$.}\label{FDUIsos}
\end{figure}

By the first isomorphism theorem, the aforementioned many-sorted  partial $\Sigma^{\boldsymbol{\mathcal{A}}}$-algebras are isomorphic.

\begin{restatable}{corollary}{CDUIsos}
\label{CDUIsos}
The many-sorted  partial $\Sigma^{\boldsymbol{\mathcal{A}}}$-algebras 
\begin{multicols}{3}
\begin{itemize}
\item[(i)] $\mathbf{Pth}^{(1,2)}_{\boldsymbol{\mathcal{A}}^{(2)}}/{\mathrm{Ker}(\mathrm{sc}^{([1],2)})}$,
\item[(ii)] $\mathbf{Pth}^{(1,2)}_{\boldsymbol{\mathcal{A}}^{(2)}}/{\mathrm{Ker}(\mathrm{tg}^{([1],2)})}$,
\item[(iii)] $[\mathbf{PT}_{\boldsymbol{\mathcal{A}}}]$.
\end{itemize}
\end{multicols}
are isomorphic.
\end{restatable}

\subsection{
\texorpdfstring
{Connecting layers $2$ and $0$}
{Connecting layers 2 and 0}
}
The structure of many-sorted $\Sigma$-algebra that we have defined on $\mathrm{Pth}_{\boldsymbol{\mathcal{A}}^{(2)}}$ has a nice interaction with the many-sorted mappings of $(0,2)$-source $(0,2)$-target and $(2,0)$-identity second-order path. 

We start by proving that the mappings  $\mathrm{sc}^{(0,2)}$ and $\mathrm{tg}^{(0,2)}$ are $\Sigma$-homomorphisms from $\mathbf{Pth}_{\boldsymbol{\mathcal{A}}^{(2)}}^{(0,2)}$ to $\mathbf{T}_{\Sigma}(X)$.

\begin{restatable}{proposition}{PDZHom}
\label{PDZHom} The mappings $\mathrm{sc}^{(0,2)}$ and $\mathrm{tg}^{(0,2)}$ are $\Sigma$-homomorphisms from $\mathbf{Pth}^{(0,2)}_{\boldsymbol{\mathcal{A}}^{(2)}}$ to $\mathbf{T}_{\Sigma}(X)$.
\end{restatable}
\begin{proof}
Let us recall from Definition~\ref{DDScTgZ} that the mapping $\mathrm{sc}^{(0,2)}$ is defined by the composition
\allowdisplaybreaks
\begin{align*}
\mathrm{sc}^{(0,2)}&=
\mathrm{sc}^{(0,[1])}\circ\mathrm{ip}^{([1],X)@}\circ\mathrm{sc}^{([1],2)}.
\end{align*}

Let us recall that $\mathrm{sc}^{([1],2)}$ is a $\Sigma^{\boldsymbol{\mathcal{A}}}$-homomorphism by Proposition~\ref{PDUCatHom}, in particular since $\Sigma\subseteq\Sigma^{\boldsymbol{\mathcal{A}}}$, we have that $\mathrm{sc}^{([1],2)}$ is a $\Sigma$-homomorphism. The free completion of the identity path mapping is also a $\Sigma$-homomorphism for the same argument. Moreover, $\mathrm{sc}^{(0,[1])}$ is a $\Sigma$-homomorphism by Proposition~\ref{PCHDZ}. All in all, we can affirm that $\mathrm{sc}^{(0,2)}$ is a $\Sigma$-homomorphism because it is a composition of $\Sigma$-homomorphisms.

Also, since we have that
$$\mathrm{tg}^{(0,2)}=
\mathrm{tg}^{(0,[1])}\circ\mathrm{ip}^{([1],X)@}\circ\mathrm{tg}^{([1],2)}$$ a similar argument applies to the $(0,2)$-target mapping.
\end{proof}

Also $\mathrm{ip}^{(2,0)\sharp}$, the $(2,0)$-identity path mapping, is a $\Sigma$-homomorphism from $\mathbf{T}_{\Sigma}(X)$ to $\mathbf{Pth}^{(0,2)}_{\boldsymbol{\mathcal{A}}^{(2)}}$.

\begin{restatable}{proposition}{PDZHomIp}
\label{PDZHomIp}
The mapping $\mathrm{ip}^{(2,0)\sharp}$ is a $\Sigma$-homomorphism from $\mathbf{T}_{\Sigma}(X)$ to $\mathbf{Pth}^{(0,2)}_{\boldsymbol{\mathcal{A}}^{(2)}}$.
\end{restatable}
\begin{proof}
Let us recall from Definition~\ref{DDPth} that $\mathrm{ip}^{(2,X)}$ is the $S$-sorted mapping from $X$ to $\mathrm{Pth}_{\boldsymbol{\mathcal{A}}^{(2)}}$ that, for every $s\in S$, sends $x\in X_{s}$ to $([x]_{s},\lambda,\lambda)$,  the $(2,0)$-identity second-order path on $[x]_{s}$.  

Since $\mathbf{Pth}^{(0,2)}_{\boldsymbol{\mathcal{A}}^{(2)}}$ is a many-sorted $\Sigma$-algebra, by the universal property of the free $\Sigma$-algebra $\mathbf{T}_{\Sigma}(X)$, there exists a unique $\Sigma$-homomorphism $\mathrm{ip}^{(2,0)\sharp}$ from $\mathbf{T}_{\Sigma}(X)$ to $\mathbf{Pth}^{(0,2)}_{\boldsymbol{\mathcal{A}}^{(2)}}$ such that
$
\mathrm{ip}^{(2,0)\sharp}\circ\eta^{(0,X)}=\mathrm{ip}^{(2,X)}, 
$
i.e., such that the diagram in Figure~\ref{FDZHomIp} commutes. 
Let us note that the $\Sigma$-homomorphism $\mathrm{ip}^{(2,0)\sharp}$ sends, for every sort $s\in S$, a  term $P\in\mathrm{T}_{\Sigma}(X)_{s}$ to $([\eta^{(1,0)\sharp}_{s}(P)]_{s},\lambda,\lambda)$, the $(2,0)$-identity path on the term $P$. 
Moreover, this definition of $\mathrm{ip}^{(2,0)\sharp}$ coincides with the presentation of $\mathrm{ip}^{(2,0)\sharp}$ introduced in Proposition~\ref{DDScTgZ}. We claim that 
$$
\mathrm{ip}^{(2,0)\sharp}=\mathrm{ip}^{(2,[1])\sharp}\circ\mathrm{CH}^{[1]}\circ\mathrm{ip}^{([1],0)\sharp}.
$$
Notice that $\mathrm{ip}^{(2,[1])\sharp}\circ\mathrm{CH}^{[1]}\circ\mathrm{ip}^{([1],0)\sharp}$ is a composition of $\Sigma$-homomorphisms in virtue of Propositions~\ref{PDUIpCatHom}, Claim~\ref{CIsoCH}, and Proposition~\ref{PCHDZIp}.

Moreover, the following chain of equalities holds
\begin{align*}
\mathrm{ip}^{(2,[1])\sharp}\circ\mathrm{CH}^{[1]}\circ\mathrm{ip}^{([1],0)\sharp}
\circ\eta^{(0,X)}
&=
\mathrm{ip}^{(2,[1])\sharp}\circ\mathrm{CH}^{[1]}
\circ\mathrm{ip}^{([1],X)}
\\&=
\mathrm{ip}^{(2,[1])\sharp}\circ\eta^{([1],X)}
\\&=
\mathrm{ip}^{(2,X)}.
\end{align*}
All of the above equalities follow from Proposition~\ref{PDBasicEq}.

So, considering the foregoing, we can affirm that 
$$
\mathrm{ip}^{(2,0)\sharp}=\mathrm{ip}^{(2,[1])\sharp}\circ\mathrm{CH}^{[1]}\circ\mathrm{ip}^{([1],0)\sharp}.$$

This finishes the proof.
\end{proof}

\begin{figure}
\begin{center}
\begin{tikzpicture}
[ACliment/.style={-{To [angle'=45, length=5.75pt, width=4pt, round]}},scale=1]
\node[] (xoq) at (0,0) [] {$X$};
\node[] (txoq) at (6,0) [] {$\mathbf{T}_{\Sigma}(X)$};
\node[] (p) at (6,-3) [] {$
\mathbf{Pth}^{(0,2)}_{\boldsymbol{\mathcal{A}}^{(2)}}
$};
\draw[ACliment]  (xoq) to node [above]
{$\eta^{(0,X)}$} (txoq);
\draw[ACliment, bend right=10]  (xoq) to node [below left] {$\mathrm{ip}^{(2,X)}$} (p);

\node[] (B0) at (6,-1.5)  [] {};
\draw[ACliment]  ($(B0)+(0,1.2)$) to node [above, fill=white] {
$\textstyle \mathrm{ip}^{(2,0)\sharp}$
} ($(B0)+(0,-1.2)$);
\draw[ACliment, bend right]  ($(B0)+(.3,-1.2)$) to node [ below, fill=white] {
$\textstyle \mathrm{tg}^{(0,2)}$
} ($(B0)+(.3,1.2)$);
\draw[ACliment, bend left]  ($(B0)+(-.3,-1.2)$) to node [below, fill=white] {
$\textstyle \mathrm{sc}^{(0,2)}$
} ($(B0)+(-.3,1.2)$);

\end{tikzpicture}
\end{center}
\caption{$\Sigma$-homomorphisms relative to $X$ at layers 0 \& 2.}
\label{FDZHomIp}
\end{figure}

From the previous results, we infer that the free $\Sigma$-algebra associated to $X$ is isomorphic to the $\Sigma$-subalgebra of the $(2,0)$-identity second-order paths.

\begin{corollary}\label{FDZIpIso} The mapping $\mathrm{ip}^{(2,0)\sharp}$ is a $\Sigma$-isomorphism from $\mathbf{T}_{\Sigma}(X)$ to $\mathrm{ip}^{(2,0)\sharp}[\mathbf{T}_{\Sigma}(X)]$.
\end{corollary}
\begin{proof}
By Proposition~\ref{PDZHomIp} $\mathrm{ip}^{(2,0)\sharp}$ is a $\Sigma$-homomorphism from $\mathbf{T}_{\Sigma}(X)$ to $\mathrm{ip}^{(2,0)\sharp}[\mathbf{T}_{\Sigma}(X)]$. The surjectivity of $\mathrm{ip}^{(2,0)\sharp}$ follows from the description of the codomain and its injectivity follows from Proposition~\ref{PDBasicEqZ}, since $\mathrm{sc}^{(0,2)}\circ\mathrm{ip}^{(2,0)\sharp}=\mathrm{id}^{\mathrm{T}_{\Sigma}(X)}$ (analogously, with the $(0,2)$-target).
\end{proof}

\begin{proposition}\label{PDZSection}
The following equations between $\Sigma$-homomorphisms hold
\begin{multicols}{2}
\begin{itemize}
\item[(i)] $\mathrm{sc}^{(0,2)}\circ \mathrm{ip}^{(2,0)\sharp}
=\mathrm{id}^{\mathbf{T}_{\Sigma}(X)}$,
\item[(ii)] $\mathrm{tg}^{(0,2)}\circ \mathrm{ip}^{(2,0)\sharp}
=\mathrm{id}^{\mathbf{T}_{\Sigma}(X)}$.
\end{itemize}
\end{multicols}
\end{proposition}
\begin{proof}
The following chain of equalities holds
\begin{align*}
\mathrm{sc}^{(0,2)}\circ\mathrm{ip}^{(2,0)\sharp}
&=
\mathrm{sc}^{(0,[1])}\circ\mathrm{ip}^{([1],X)@}\circ\mathrm{sc}^{([1],2)}
\circ
\mathrm{ip}^{(2,[1])\sharp}\circ\mathrm{CH}^{[1]}\circ\mathrm{ip}^{([1],0)\sharp}
\tag{1}
\\&=
\mathrm{sc}^{(0,[1])}\circ\mathrm{ip}^{([1],X)@}\circ\mathrm{CH}^{[1]}\circ\mathrm{ip}^{([1],0)\sharp}
\tag{2}
\\&=
\mathrm{sc}^{(0,[1])}\circ\mathrm{ip}^{([1],0)\sharp}
\tag{3}
\\&=
\mathrm{id}^{\mathbf{T}_{\Sigma}(X)}.
\tag{4}
\end{align*}

The first equality unpacks the definition of the mappings $\mathrm{sc}^{(0,2)}$ and $\mathrm{ip}^{(2,0)\sharp}$ introduced in Definition~\ref{DDScTgZ}; the second equality follows from associativity and from Proposition~\ref{PDBasicEq}; the third equality follows from associativity and from Theorem~\ref{TIso}; finally, the last equality follows from Proposition~\ref{PCHBasicEq}.

A similar argument can be used to obtain that 
$$
\mathrm{tg}^{(0,2)}\circ \mathrm{ip}^{(2,0)\sharp}
=\mathrm{id}^{\mathbf{T}_{\Sigma}(X)}.
$$

This finishes the proof.
\end{proof}

Thus, $\mathbf{T}_{\Sigma}(X)$ is a retract of $\mathbf{Pth}^{(0,2)}_{\boldsymbol{\mathcal{A}}^{(2)}}$, the
$\Sigma$-homomorphisms $\mathrm{sc}^{(0,2)}$ and $\mathrm{tg}^{(0,2)}$ are retractions of $\mathbf{Pth}^{(0,2)}_{\boldsymbol{\mathcal{A}}^{(2)}}$ onto $\mathbf{T}_{\Sigma}(X)$, and the $\Sigma$-homomorphism $\mathrm{ip}^{(2,0)\sharp}$ is a section of both $\mathrm{sc}^{(0,2)}$ and $\mathrm{tg}^{(0,2)}$.

\begin{corollary}\label{CDZScTg}
Let $s$ be a sort in $S$ and $\mathfrak{P}^{(2)}$ a second-order path in $\mathrm{Pth}_{\boldsymbol{\mathcal{A}}^{(2)},s}$. Then we have that
\begin{align*}
\left(\mathfrak{P}^{(2)},\mathrm{ip}^{(2,0)\sharp}_{s}
\left(\mathrm{sc}^{(0,2)}_{s}
\left(\mathfrak{P}^{(2)}
\right)\right)\right)
&\in\mathrm{Ker}\left(\mathrm{sc}^{(0,2)}\right)_{s},\\
\left(\mathfrak{P}^{(2)},\mathrm{ip}^{(2,0)\sharp}_{s}
\left(\mathrm{tg}^{(0,2)}_{s}
\left(\mathfrak{P}^{(2)}
\right)\right)\right)
&\in\mathrm{Ker}\left(\mathrm{tg}^{(0,2)}\right)_{s}.
\end{align*}
\end{corollary}

\begin{remark}
From Proposition~\ref{PDZHom} it follows that the kernel of $\mathrm{sc}^{(0,2)}$ is a $\Sigma$-congruence on $\mathbf{Pth}^{(0,2)}_{\boldsymbol{\mathcal{A}}^{(2)}}$. Therefore, we can consider the standard (epi, mono)-factoriza\-tion of $\mathrm{sc}^{(0,2)}$ given by
\begin{enumerate}
\item the projection $\mathrm{pr}^{\mathrm{Ker}(\mathrm{sc}^{(0,2)})}$ from $\mathbf{Pth}^{(0,2)}_{\boldsymbol{\mathcal{A}}^{(2)}}$ to $\mathbf{Pth}^{(0,2)}_{\boldsymbol{\mathcal{A}}^{(2)}}/{\mathrm{Ker}(\mathrm{sc^{(0,2)}})}$ that, for every $s\in S$, assigns to a second-order path $\mathfrak{P}^{(2)}$ in $\mathrm{Pth}^{(0,2)}_{\boldsymbol{\mathcal{A}}^{(2)},s}$ its equivalence class $[\mathfrak{P}^{(2)}]_{\mathrm{Ker}(\mathrm{sc}^{(0,2)})_{s}}$, and
\item the embedding $\mathrm{sc}^{(0,2)\mathrm{m}}$ of $\mathbf{Pth}^{(0,2)}_{\boldsymbol{\mathcal{A}}^{(2)}}/{\mathrm{Ker}(\mathrm{sc^{(0,2)}})}$ into $\mathbf{T}_{\Sigma}(X)$ that, for every $s\in S$, assigns to an equivalence class $[\mathfrak{P}^{(2)}]_{\mathrm{Ker}(\mathrm{sc}^{(0,2)})_{s}}$ in $\mathrm{Pth}^{(0,2)}_{\boldsymbol{\mathcal{A}}^{(2)},s}/{\mathrm{Ker}(\mathrm{sc}^{(0,2)})_{s}}$, with $\mathfrak{P}^{(2)}\in\mathrm{Pth}_{\boldsymbol{\mathcal{A}}^{(2)},s}$, the term $\mathrm{sc}^{(0,2)}_{s}(\mathfrak{P}^{(2)})$, i.e., the value of the mapping $\mathrm{sc}^{(0,2)}_{s}$ at any equivalence class representative. We will refer to the mapping $\mathrm{sc}^{(0,2)\mathrm{m}}$ as the \emph{monomorphic part of the factorization of $\mathrm{sc}^{(0,2)}$}.
\end{enumerate}

Let us recall that $(\mathrm{pr}^{\mathrm{Ker}(\mathrm{sc}^{(0,2)})}, \mathrm{sc}^{(0,2)\mathrm{m}})$ is the unique pair of $\Sigma$-homomorphisms such that $\mathrm{sc}^{(0,2)\mathrm{m}}\circ\mathrm{pr}^{\mathrm{Ker}(\mathrm{sc^{(0,2)}})}=\mathrm{sc}^{(0,2)}$, i.e., such that the upper half of the diagram in Figure~\ref{FDZIsos} commutes.

On the other hand, the kernel of $\mathrm{tg}^{(0,2)}$ is also a $\Sigma$-congruence on $\mathbf{Pth}^{(0,2)}_{\boldsymbol{\mathcal{A}}^{(2)}}$. We can consider the standard (epi, mono)-factorization of $\mathrm{tg}^{(0,2)}$ given by
\begin{enumerate}
\item the projection $\mathrm{pr}^{\mathrm{Ker}(\mathrm{tg}^{(0,2)})}$ from $\mathbf{Pth}^{(0,2)}_{\boldsymbol{\mathcal{A}}^{(2)}}$ to $\mathbf{Pth}^{(0,2)}_{\boldsymbol{\mathcal{A}}^{(2)}}/{\mathrm{Ker}(\mathrm{tg}^{(0,2)})}$ that, for every $s\in S$, assigns to a second-order path $\mathfrak{P}^{(2)}$ in $\mathrm{Pth}_{\boldsymbol{\mathcal{A}}^{(2)},s}$ its equivalence class $[\mathfrak{P}^{(2)}]_{\mathrm{Ker}(\mathrm{tg}^{(0,2)})_{s}}$, and
\item the embedding $\mathrm{tg}^{(0,2)\mathrm{m}}$ from $\mathbf{Pth}^{(0,2)}_{\boldsymbol{\mathcal{A}}^{(2)}}/{\mathrm{Ker}(\mathrm{tg}^{(0,2)})}$ to $\mathbf{T}_{\Sigma}(X)$ that, for every $s\in S$, assigns to an equivalence class $[\mathfrak{P}^{(2)}]_{\mathrm{Ker}(\mathrm{tg}^{(0,2)})_{s}}$ in $\mathrm{Pth}_{\boldsymbol{\mathcal{A}}^{(2)},s}/{\mathrm{Ker}(\mathrm{tg}^{(0,2)})_{s}}$ and $\mathfrak{P}^{(2)}\in\mathrm{Pth}_{\boldsymbol{\mathcal{A}}^{(2)},s}$ the term $\mathrm{tg}^{(0,2)}_{s}(\mathfrak{P}^{(2)})$, i.e., the value of the mapping $\mathrm{tg}^{(0,2)}_{s}$ at any equivalence class representative. We will refer to the mapping $\mathrm{tg}^{(0,2)\mathrm{m}}$ as the \emph{monomorphic part of the factorization of $\mathrm{tg}^{(0,2)}$}.
\end{enumerate}

Let us recall that $(\mathrm{pr}^{\mathrm{Ker}(\mathrm{tg}^{(0,2)})}, \mathrm{tg}^{(0,2)\mathrm{m}})$ is the unique pair of $\Sigma$-homomorphisms such that $\mathrm{tg}^{(0,2)\mathrm{m}}\circ\mathrm{pr}^{\mathrm{Ker}(\mathrm{tg}^{(0,2)})}=\mathrm{tg}^{(0,2)}$, i.e., such that the lower half of the diagram in Figure~\ref{FDZIsos} commutes.
\end{remark}

\begin{figure}
\begin{center}
\begin{tikzpicture}
[ACliment/.style={-{To [angle'=45, length=5.75pt, width=4pt, round]}}]
\node[] (txoq) at (-1,0) [] {$\mathbf{T}_{\Sigma}(X)$};
\node[] (p) at (4,3.5) [] {$\mathbf{Pth}^{(0,2)}_{\boldsymbol{\mathcal{A}}^{(2)}}$};
\node[] (p2) at (4,-3.5) [] {$\mathbf{Pth}^{(0,2)}_{\boldsymbol{\mathcal{A}}^{(2)}}$};
\node[] (pkernel) at (4,1.5) [] {$\mathbf{Pth}^{(0,2)}_{\boldsymbol{\mathcal{A}}^{(2)}}/{\mathrm{Ker}(\mathrm{sc}^{(0,2)})}$};
\node[] (pkernel2) at (4,-1.5) [] {$\mathbf{Pth}^{(0,2)}_{\boldsymbol{\mathcal{A}}^{(2)}}/{\mathrm{Ker}(\mathrm{tg}^{(0,2)})}$};
\node[] (txoq2) at (9,0) [] {$\mathbf{T}_{\Sigma}(X)$};
\draw[ACliment, bend left=15]  (txoq) to node [above left]
{$\mathrm{ip}^{(2,0)\sharp}$} (p);
\draw[ACliment]  (txoq) to node [sloped, above] {}
(pkernel);
\draw[ACliment]  (p) to node [midway, fill=white]
{$\mathrm{pr}^{\mathrm{Ker}(\mathrm{sc}^{(0,2)})}$} (pkernel);
\draw[ACliment]  (pkernel) to node [sloped, above]
{$\mathrm{sc}^{(0,2)\mathrm{m}}$} (txoq2);
\draw[ACliment, bend left=15]  (p) to node [above right]
{$\mathrm{sc}^{(0,2)}$} (txoq2);
\draw[ACliment]  (txoq) to node  [midway, fill=white]
{$\mathrm{id}^{\mathbf{T}_{\Sigma}(X)}$} (txoq2);
\draw[ACliment, bend right=15]  (txoq) to node [below left]
{$\mathrm{ip}^{(2,0)\sharp}$} (p2);
\draw[ACliment]  (txoq) to node [sloped, below] {}
(pkernel2);
\draw[ACliment]  (p2) to node  [midway, fill=white]
{$\mathrm{pr}^{\mathrm{Ker}(\mathrm{tg}^{(0,2)})}$} (pkernel2);
\draw[ACliment]  (pkernel2) to node [sloped, below]
{$\mathrm{tg}^{(0,2)\mathrm{m}}$} (txoq2);
\draw[ACliment, bend right=15]  (p2) to node [below right]
{$\mathrm{tg}^{(0,2)}$} (txoq2);
\end{tikzpicture}
\end{center}
\caption{$\Sigma$-homomorphisms and $\Sigma$-algebras at layers $2$ \& $0$.}\label{FDZIsos}
\end{figure}

By the first isomorphism theorem, the aforementioned many-sorted  $\Sigma$-algebras are isomorphic.

\begin{restatable}{corollary}{CDZIsos}
\label{CDZIsos}
The many-sorted  $\Sigma$-algebras 
\begin{multicols}{3}
\begin{itemize}
\item[(i)] $\mathbf{Pth}^{(0,2)}_{\boldsymbol{\mathcal{A}}^{(2)}}/{\mathrm{Ker}(\mathrm{sc}^{(0,2)})}$,
\item[(ii)] $\mathbf{Pth}^{(0,2)}_{\boldsymbol{\mathcal{A}}^{(2)}}/{\mathrm{Ker}(\mathrm{tg}^{(0,2)})}$,
\item[(iii)] $\mathbf{T}_{\Sigma}(X)$.
\end{itemize}
\end{multicols}
are isomorphic.
\end{restatable}			
\chapter{
\texorpdfstring
{An Artinian order on $\coprod \mathrm{Pth}_{\boldsymbol{\mathcal{A}}^{(2)}}$}
{An Artinian order on second-order paths}
}\label{S2C}

In this chapter we define on $\coprod \mathrm{Pth}_{\boldsymbol{\mathcal{A}}^{(2)}}$, the coproduct of $\mathrm{Pth}_{\boldsymbol{\mathcal{A}}^{(2)}}$---formed by all labelled second-order paths $(\mathfrak{P}^{(2)},s)$, with $s\in S$ and $\mathfrak{P}^{(2)}\in \mathrm{Pth}_{\boldsymbol{\mathcal{A}}^{(2)},s}$---, an Artinian order (see Definition~\ref{DPosArt}), which will allow us to justify both proofs by Artinian induction and definitions by Artinian recursion, which we will use in subsequent chapters. Moreover, we specify the minimals of the corresponding Artinian ordered set, prove that $\mathrm{ip}^{(2,[1])\sharp}[[\mathrm{PT}_{\boldsymbol{\mathcal{A}}}]]$ is a lower set and that several mappings to $\coprod \mathrm{Pth}_{\boldsymbol{\mathcal{A}}^{(2)}}$ are order embeddings.



\begin{restatable}{definition}{DDOrd}
\label{DDOrd}
\index{partial order!second-order!$\leq_{\mathbf{Pth}_{\boldsymbol{\mathcal{A}}^{(2)}}}$}
We let $\prec_{\mathbf{Pth}_{\boldsymbol{\mathcal{A}}^{(2)}}}$ denote the binary relation on $\coprod \mathrm{Pth}_{\boldsymbol{\mathcal{A}}^{(2)}}$ consisting of the ordered pairs $((\mathfrak{Q}^{(2)},t),(\mathfrak{P}^{(2)},s))\in(\coprod \mathrm{Pth}_{\boldsymbol{\mathcal{A}}^{(2)}})^{2}$ for which one of the following conditions holds
\begin{enumerate}
\item $\mathfrak{P}^{(2)}$ and $\mathfrak{Q}^{(2)}$ are $(2,[1])$-identity second-order paths of the form 
\begin{align*}
\mathfrak{P}^{(2)}&=\mathrm{ip}^{(2,[1])\sharp}_{s}\left([P]_{s}\right),
&
\mathfrak{Q}^{(2)}&=\mathrm{ip}^{(2,[1])\sharp}_{t}\left([Q]_{t}\right),
\end{align*}
for some path term classes $[P]_{s}\in [\mathrm{PT}_{\boldsymbol{\mathcal{A}}}]_{s}$ and $[Q]_{t}\in [\mathrm{PT}_{\boldsymbol{\mathcal{A}}}]_{t}$, and the following inequality holds
$$
\left([Q]_{t},t
\right)
<_{[\mathbf{PT}_{\boldsymbol{\mathcal{A}}}]}
\left([P]_{s},s
\right),
$$
where $\leq_{[\mathbf{PT}_{\boldsymbol{\mathcal{A}}}]}$ is the Artinian partial order on $\coprod[\mathrm{PT}_{\boldsymbol{\mathcal{A}}}]$ introduced in Definition~\ref{DPTQOrd}.
\item $\mathfrak{P}^{(2)}$ is a second-order path of length strictly greater than one containing at least one   second-order echelon, and if its first   second-order echelon occurs at position $i\in\bb{\mathfrak{P}^{(2)}}$, then
\subitem(2.1) if $i=0$, then $\mathfrak{Q}^{(2)}$ is equal to $\mathfrak{P}^{(2),0,0}$ or $\mathfrak{P}^{(2),1,\bb{\mathfrak{P}^{(2)}}-1}$,
\subitem(2.2) if $i>0$, then $\mathfrak{Q}^{(2)}$ is equal to $\mathfrak{P}^{(2),0,i-1}$ or $\mathfrak{P}^{(2),i,\bb{\mathfrak{P}^{(2)}}-1}$,
\item $\mathfrak{P}^{(2)}$ is an echelonless second-order path, then
\subitem(3.1) if $\mathfrak{P}^{(2)}$ is not head-constant, then let $i\in\bb{\mathfrak{P}^{(2)}}$ be the  maximum index for which $\mathfrak{P}^{(2),0,i}$ is a head-constant second-order path, then $\mathfrak{Q}^{(2)}$ is equal to $\mathfrak{P}^{(2),0,i}$ or $\mathfrak{P}^{(2),i+1,\bb{\mathfrak{P}^{(2)}}-1}$,
\subitem(3.2) if $\mathfrak{P}^{(2)}$ is head-constant but not coherent, then let $i\in\bb{\mathfrak{P}^{(2)}}$ be the  maximum index for which $\mathfrak{P}^{(2),0,i}$ is a coherent second-order path, then $\mathfrak{Q}^{(2)}$ is equal to $\mathfrak{P}^{(2),0,i}$ or $\mathfrak{P}^{(2),i+1,\bb{\mathfrak{P}^{(2)}}-1}$,
\subitem(3.3) if $\mathfrak{P}^{(2)}$ is head-constant and coherent then $\mathfrak{Q}^{(2)}$ is one of the second-order paths we can extract from $\mathfrak{P}^{(2)}$ in virtue of Lemma~\ref{LDPthExtract}.
\end{enumerate}
The reader is advised to consult the flowchart from Figure~\ref{FDFlow} for deciding on whether, given two pairs $(\mathfrak{Q}^{(2)},t)$ and $(\mathfrak{P}^{(2)},s)$ in $\coprod\mathrm{Pth}_{\boldsymbol{\mathcal{A}}^{(2)}}$, the statement $(\mathfrak{Q}^{(2)},t)\prec_{\mathbf{Pth}_{\boldsymbol{\mathcal{A}}^{(2)}}}(\mathfrak{P}^{(2)},s)$ is fulfilled or not fulfilled.

We will denote by $\leq_{\mathbf{Pth}_{\boldsymbol{\mathcal{A}}^{(2)}}}$ the reflexive and transitive closure of $\prec_{\mathbf{Pth}_{\boldsymbol{\mathcal{A}}^{(2)}}}$, i.e., the preorder on $\coprod \mathrm{Pth}_{\boldsymbol{\mathcal{A}}^{(2)}}$  generated by $\prec_{\mathbf{Pth}_{\boldsymbol{\mathcal{A}}^{(2)}}}$, and by $<_{\mathbf{Pth}_{\boldsymbol{\mathcal{A}}^{(2)}}}$ the transitive closure of $\prec_{\mathbf{Pth}_{\boldsymbol{\mathcal{A}}^{(2)}}}$.
\end{restatable}

\begin{figure}
\begin{tikzpicture}[scale=1.1, baseline=0cm]
\tikzstyle{end} =  [regular polygon,regular polygon sides=5, draw, shape border rotate=180,  minimum width=.8cm, minimum height=.8cm, text centered, draw=black, fill=purple!20!white]
\tikzstyle{start} =  [regular polygon,regular polygon sides=5, draw, minimum width=.8cm, minimum height=.8cm, text centered, draw=black, fill=purple!20!white]
\tikzstyle{true} = [rectangle, rounded corners, minimum width=.7cm, minimum height=.7cm, text centered, draw=black, fill=green!20!white]
\tikzstyle{false} = [rectangle, rounded corners, minimum width=.7cm, minimum height=.7cm, text centered, draw=black, fill=red!20!white]
\tikzstyle{io} = [trapezium, trapezium left angle=70, trapezium right angle=110, minimum width=.55cm, minimum height=.55cm, text centered, draw=black, fill=purple!20!white]
\tikzstyle{process} = [rectangle, minimum width=.7cm, minimum height=.7cm, text centered, draw=black, fill=yellow!20!white]
\tikzstyle{decision} = [diamond, minimum width=.8cm, minimum height=.8cm, text centered, draw=black, fill=blue!20!white]
\tikzstyle{arrow} = [thick,->,>=stealth]

\node (A1) at (2,-1) [start] {$\scriptstyle 1$};
\node (A2) at (2,-2) [decision] {$\scriptstyle 2$};
\node (A3) at (-1,-3) [decision] {$\scriptstyle 3$};
\node (A4) at (4,-3) [decision] {$\scriptstyle 4$};
\node (A5) at (-2,-4) [process] {$\scriptstyle 5$};
\node (F0) at (-1,-4) [false] {$\scriptstyle \mathrm{F}$};
\node (A6) at (2,-4) [decision] {$\scriptstyle 6$};
\node (F1) at (2,-5) [false] {$\scriptstyle \mathrm{F}$};
\node (A7) at (6,-4) [decision] {$\scriptstyle 7$};
\node (A8) at (-2,-5) [end] {$\scriptstyle 8$};
\node (A9) at (0,-5) [process] {$\scriptstyle 9$};
\node (A10) at (4,-5) [decision] {$\scriptstyle 10$};
\node (A11) at (7,-5) [process] {$\scriptstyle 11$};
\node (A12) at (0,-6) [decision] {$\scriptstyle 12$};
\node (A13) at (3,-6) [decision] {$\scriptstyle 13$};
\node (A14) at (5,-6) [process] {$\scriptstyle 14$};
\node (A15) at (7,-6) [decision] {$\scriptstyle 15$};
\node (A16) at (-1,-7) [decision] {$\scriptstyle 16$};
\node (V1) at (-3,-6) [true] {$\scriptstyle \mathrm{T}$};
\node (F2) at (-2,-6) [false] {$\scriptstyle \mathrm{F}$};
\node (A17) at (1,-7) [decision] {$\scriptstyle 17$};
\node (V2) at (-2,-8) [true] {$\scriptstyle \mathrm{T}$};
\node (F3) at (-1,-8) [false] {$\scriptstyle \mathrm{F}$};
\node (V3) at (0,-8) [true] {$\scriptstyle \mathrm{T}$};
\node (F4) at (1,-8) [false] {$\scriptstyle \mathrm{F}$};
\node (V4) at (2,-7) [true] {$\scriptstyle \mathrm{T}$};
\node (F5) at (3,-7) [false] {$\scriptstyle \mathrm{F}$};
\node (A18) at (5,-7) [decision] {$\scriptstyle 18$};
\node (V5) at (6,-7) [true] {$\scriptstyle \mathrm{T}$};
\node (F6) at (7,-7) [false] {$\scriptstyle \mathrm{F}$};
\node (V6) at (4,-8) [true] {$\scriptstyle \mathrm{T}$};
\node (F7) at (5,-8) [false] {$\scriptstyle \mathrm{F}$};

\draw [arrow] (A1) -- (A2);
\draw [arrow] (A2) -| node[anchor=south west] {yes} (A3);
\draw [arrow] (A2) -| node[anchor=south east] {no} (A4);
\draw [arrow] (A3) -| node[anchor=south west] {yes} (A5);
\draw [arrow] (A3) -- node[anchor=west] {no} (F0);
\draw [arrow] (A4) -| node[anchor=south west] {yes} (A6);
\draw [arrow] (A4) -| node[anchor=south east] {no} (A7);
\draw [arrow] (A6) -| node[anchor=south west] {yes} (A9);
\draw [arrow] (A6) -- node[anchor=west] {no} (F1);
\draw [arrow] (A5) -- (A8);
\draw [arrow] (A9) -- (A12);
\draw [arrow] (A7) -| node[anchor=south west] {yes} (A10);
\draw [arrow] (A7) -| node[anchor=south east] {no} (A11);
\draw [arrow] (A8) -| node[anchor=south west] {yes} (V1);
\draw [arrow] (A8) -- node[anchor=west] {no} (F2);
\draw [arrow] (A12) -| node[anchor=south west] {yes} (A16);
\draw [arrow] (A12) -| node[anchor=south east] {no} (A17);
\draw [arrow] (A10) -| node[anchor=south west] {yes} (A13);
\draw [arrow] (A10) -| node[anchor=south east] {no} (A14);
\draw [arrow] (A11) -- (A15);
\draw [arrow] (A16) -| node[anchor=south west] {yes} (V2);
\draw [arrow] (A16) -- node[anchor=west] {no} (F3);
\draw [arrow] (A17) -| node[anchor=south west] {yes} (V3);
\draw [arrow] (A17) -- node[anchor=west] {no} (F4);
\draw [arrow] (A13) -| node[anchor=south west] {yes} (V4);
\draw [arrow] (A13) -- node[anchor=west] {no} (F5);
\draw [arrow] (A14) -- (A18);
\draw [arrow] (A15) -| node[anchor=south west] {yes} (V5);
\draw [arrow] (A15) -- node[anchor=west] {no} (F6);
\draw [arrow] (A18) -| node[anchor=south west] {yes} (V6);
\draw [arrow] (A18) -- node[anchor=west] {no} (F7);
\end{tikzpicture}
\begin{center}
\begin{tabu} to 1.2\textwidth {  r  X[l]  }
\\
1.&Let $(\mathfrak{Q}^{(2)},t),(\mathfrak{P}^{(2)},s)$ be second-order paths in $\coprod\mathrm{Pth}_{\boldsymbol{\mathcal{A}}^{(2)}}$\\
2.&Is $\mathfrak{P}^{(2)}$ a $(2,[1])$-identity second-order path?\\
3.&Is $\mathfrak{Q}^{(2)}$ a $(2,[1])$-identity second-order path?\\
4.&Does $\mathfrak{P}^{(2)}$ contain at least one   second-order echelon?\\
5.&Let $\mathfrak{Q}^{(2)}=\mathrm{ip}^{(2,[1])\sharp}_{t}([Q]_{t})$ and
$\mathfrak{P}^{(2)}=\mathrm{ip}^{(2,[1])\sharp}_{t}([P]_{t})$\\
6.&Is $\bb{\mathfrak{P}^{(2)}}>1$?\\
7.&Is $\mathfrak{P}^{(2)}$ a head-constant?\\
8.&Is $([Q]_{t},t)\leq_{[\mathbf{PT}_{\boldsymbol{\mathcal{A}}}]}([P]_{s},s)$?\\
9.&Let $i\in\bb{\mathfrak{P}^{(2)}}$ be the first index for which $\mathfrak{P}^{(2),i,i}$ is a   second-order echelon\\
10.&Is $\mathfrak{P}^{(2)}$ coherent?\\
11.&Let $i\in\bb{\mathfrak{P}^{(2)}}$ be the maximum index for which $\mathfrak{P}^{(2),0,i}$ is head-constant\\
12.&Is $i=0$?\\
13.&Is $\mathfrak{Q}^{(2)}$ one of the second-order paths we can extract from $\mathfrak{P}^{(2)}$?\\
14.&Let $i\in\bb{\mathfrak{P}^{(2)}}$ be the maximum index for which $\mathfrak{P}^{(2),0,i}$ is coherent\\
15.&Is $\mathfrak{Q}^{(2)}$ equal to $\mathfrak{P}^{(2),0,i}$ or $\mathfrak{P}^{(2),i+1,\bb{\mathfrak{P}^{(2)}}-1}$?\\
16.&Is $\mathfrak{Q}^{(2)}$ equal to $\mathfrak{P}^{(2),0,0}$ or $\mathfrak{P}^{(2),1,\bb{\mathfrak{P}^{(2)}}-1}$?\\
17.&Is $\mathfrak{Q}^{(2)}$ equal to $\mathfrak{P}^{(2),0,i-1}$ or $\mathfrak{P}^{(2),i,\bb{\mathfrak{P}^{(2)}}-1}$?\\
18.&Is $\mathfrak{Q}^{(2)}$ equal to $\mathfrak{P}^{(2),0,i}$ or $\mathfrak{P}^{(2),i+1,\bb{\mathfrak{P}^{(2)}}-1}$?\\
\\
\end{tabu}
\end{center}
\begin{tikzpicture}[scale=1.1, baseline=0cm]
\tikzstyle{true} = [rectangle, rounded corners, minimum width=.7cm, minimum height=.7cm, text centered, draw=black, fill=green!20!white]
\tikzstyle{false} = [rectangle, rounded corners, minimum width=.7cm, minimum height=.7cm, text centered, draw=black, fill=red!20!white]

\node () at (0,0) [true] {$\scriptstyle \mathrm{T}$};
\node () at (2.35,0) [] {$~(\mathfrak{Q}^{(2)},t)\prec_{\mathbf{Pth}_{\boldsymbol{\mathcal{A}}^{(2)}}}(\mathfrak{P}^{(2)},s)$};
\node () at (6,0) [false] {$\scriptstyle \mathrm{F}$};
\node () at (8.35,0) [] {$~(\mathfrak{Q}^{(2)},t)\not\prec_{\mathbf{Pth}_{\boldsymbol{\mathcal{A}}^{(2)}}}(\mathfrak{P}^{(2)},s)$};
\end{tikzpicture}
\caption{Decision flowchart for $(\mathfrak{Q}^{(2)},t)\prec_{\mathbf{Pth}_{\boldsymbol{\mathcal{A}}^{(2)}}}(\mathfrak{P}^{(2)},s)$.}\label{FDFlow}
\end{figure}

\begin{restatable}{remark}{RDIrrefl}
\label{RDIrrefl}
The binary relation $\prec_{\mathbf{Pth}_{\boldsymbol{\mathcal{A}}^{(2)}}}$ on $\coprod \mathrm{Pth}_{\boldsymbol{\mathcal{A}}^{(2)}}$ is irreflexive.
\end{restatable}

\begin{restatable}{remark}{RDOrd}
\label{RDOrd}
The preorder $\leq_{\mathbf{Pth}_{\boldsymbol{\mathcal{A}}^{(2)}}}$ on $\coprod \mathrm{Pth}_{\boldsymbol{\mathcal{A}}^{(2)}}$ is $\bigcup_{n\in\mathbb{N}}\prec_{\mathbf{Pth}_{\boldsymbol{\mathcal{A}}^{(2)}}}^{n}$, where $\prec_{\mathbf{Pth}_{\boldsymbol{\mathcal{A}}^{(2)}}}^{0}$ is the diagonal of $\coprod \mathrm{Pth}_{\boldsymbol{\mathcal{A}}^{(2)}}$ and, for $n\in\mathbb{N}$, $$\prec_{\mathbf{Pth}_{\boldsymbol{\mathcal{A}}^{(2)}}}^{n+1} = \prec_{\mathbf{Pth}_{\boldsymbol{\mathcal{A}}^{(2)}}}^{n}\circ \prec_{\mathbf{Pth}_{\boldsymbol{\mathcal{A}}^{(2)}}}.$$ 

Thus, for every $((\mathfrak{Q}^{(2)},t),(\mathfrak{P}^{(2)},s))\in (\coprod \mathrm{Pth}_{\boldsymbol{\mathcal{A}}^{(2)}})^{2}$, $(\mathfrak{Q}^{(2)},t)\leq_{\mathbf{Pth}_{\boldsymbol{\mathcal{A}}^{(2)}}}(\mathfrak{P}^{(2)},s)$ if and only if $s = t$ and $\mathfrak{Q}^{(2)} = \mathfrak{P}^{(2)}$ or there exists  a natural number $m\in\mathbb{N}-\{0\}$, a word $\mathbf{w}\in S^{\star}$ of length $\bb{\mathbf{w}}=m+1$, and a family of second-order paths $(\mathfrak{R}^{(2)}_{k})_{k\in\bb{\mathbf{w}}}$ in $\mathrm{Pth}_{\boldsymbol{\mathcal{A}}^{(2)},{\mathbf{w}}}$, such that $w_{0}=t$,  $\mathfrak{R}^{(2)}_{0}=\mathfrak{Q}^{(2)}$, $w_{m}=s$, $\mathfrak{R}^{(2)}_{m}=\mathfrak{P}^{(2)}$ and, for every $k\in m$, $(\mathfrak{R}^{(2)}_{k},w_{k})\prec_{\mathbf{Pth}_{\boldsymbol{\mathcal{A}}^{(2)}}}(\mathfrak{R}^{(2)}_{k+1},w_{k+1})$. Moreover, $<_{\mathbf{Pth}_{\boldsymbol{\mathcal{A}}^{(2)}}}$, the transitive closure of $\prec_{\mathbf{Pth}_{\boldsymbol{\mathcal{A}}^{(2)}}}$, is $\bigcup_{n\in\mathbb{N}-
\{0\}}\prec_{\mathbf{Pth}_{\boldsymbol{\mathcal{A}}^{(2)}}}^{n}$.
\end{restatable}

In the following proposition we show that $\mathrm{Min}(\coprod \mathrm{Pth}_{\boldsymbol{\mathcal{A}}^{(2)}}, \leq_{\mathbf{Pth}_{\boldsymbol{\mathcal{A}}^{(2)}}})$, the set of the minimal elements of the preordered set $(\coprod \mathrm{Pth}_{\boldsymbol{\mathcal{A}}^{(2)}}, \leq_{\mathbf{Pth}_{\boldsymbol{\mathcal{A}}^{(2)}}})$, consists of the $(2,[1])$-identity second-order paths on minimal path term classes (with respect to $\leq_{[\mathbf{PT}_{\boldsymbol{\mathcal{A}}}]}$ together with the  labelled second-order echelons in $\coprod\mathrm{Pth}_{\boldsymbol{\mathcal{A}}^{(2)}}$.

\begin{restatable}{proposition}{PDMinimal}
\label{PDMinimal}
For the preordered set $(\coprod \mathrm{Pth}_{\boldsymbol{\mathcal{A}}^{(2)}}, \leq_{\mathbf{Pth}_{\boldsymbol{\mathcal{A}}^{(2)}}})$ we have that
\begin{flushleft}
$\textstyle
\mathrm{Min}\left(\coprod \mathrm{Pth}_{\boldsymbol{\mathcal{A}}^{(2)}}, \leq_{\mathbf{Pth}_{\boldsymbol{\mathcal{A}}^{(2)}}}\right)$
\begin{multline*}
\qquad=
\textstyle
\left(\coprod\mathrm{ip}^{(2,[1])\sharp}\left[
\mathrm{Min}\left(\coprod [\mathrm{PT}_{\boldsymbol{\mathcal{A}}}], \leq_{[\mathbf{PT}_{\boldsymbol{\mathcal{A}}}]}\right)\right]\right)
\cup
\textstyle
\left(\coprod\mathrm{Ech}^{(2,\mathcal{A}^{(2)})}[\mathcal{A}^{(2)}]\right).
\end{multline*}
\end{flushleft}
\end{restatable}

%

\begin{proof}
We first prove that both labelled $(2,[1])$-identity second-order paths on minimal path term classes and labelled   second-order  echelons are minimal elements.

Let $s$ be a sort in $S$ and $[P]_{s}$ a path term class in $[\mathrm{PT}_{\boldsymbol{\mathcal{A}}}]_{s}$.  Assume that $([P]_{s},s)$ is a minimal element in $(\coprod[\mathrm{PT}_{\boldsymbol{\mathcal{A}}}], \leq_{[\mathbf{PT}_{\boldsymbol{\mathcal{A}}}]})$. Consider the $(2,[1])$-identity second-order path on $[P]_{s}$, i.e., $\mathrm{ip}^{(2,[1])\sharp}_{s}([P]_{s})$. Let us note that $\mathrm{ip}^{(2,[1])\sharp}_{s}([P]_{s})$ is a second-order path of length $0$. Therefore, for every sort $t\in S$ and every second-order path $\mathfrak{Q}^{(2)}\in\mathrm{Pth}_{\boldsymbol{\mathcal{A}}^{(2)},t}$, we have that $(\mathfrak{Q}^{(2)},t)\prec_{\mathbf{Pth}_{\boldsymbol{\mathcal{A}}^{(2)}}} (\mathrm{ip}^{(2,[1])\sharp}_{s}([P]_{s}),s)$, if and only if, $\mathfrak{Q}^{(2)}$ is a $(2,[1])$-identity second-order path of the form $\mathfrak{Q}^{(2)}=\mathrm{ip}^{(2,[1])\sharp}_{t}([Q]_{t})$, for some path term class $[Q]_{t}$ in $[\mathrm{PT}_{\boldsymbol{\mathcal{A}}}]_{t}$ and 
$$
\left([Q]_{t},t
\right)
\leq_{[\mathbf{PT}_{\boldsymbol{\mathcal{A}}}]}
\left([P]_{s},s
\right).
$$
Since $([P]_{s},s)$ is minimal in $(\coprod[\mathrm{PT}_{\boldsymbol{\mathcal{A}}}], \leq_{[\mathbf{PT}_{\boldsymbol{\mathcal{A}}}]})$, we have that $t=s$ and $[Q]_{t}=[P]_{s}$. Therefore, we have that $s=t$ and
$$
\mathfrak{Q}^{(2)}
=
\mathrm{ip}^{(2,[1])\sharp}_{t}
\left([Q]_{t}
\right)
=
\mathrm{ip}^{(2,[1])\sharp}_{s}\left(
[P]_{s}
\right).
$$

It follows that $(\mathrm{ip}^{(2,[1])\sharp}_{s}([P]_{s}),s)$ is minimal in $(\coprod \mathrm{Pth}_{\boldsymbol{\mathcal{A}}^{(2)}}, \leq_{\mathbf{Pth}_{\boldsymbol{\mathcal{A}}^{(2)}}})$.

Let $s$ be a sort in $S$ and $\mathfrak{p}^{(2)}$ a second-order rewrite rule in $\mathcal{A}^{(2)}_{s}$ of the form $\mathfrak{p}^{(2)}=([M]_{s}, [N]_{s})$. Consider the second-order echelon associated to $\mathfrak{p}^{(2)}$
$$
\xymatrix@C=90pt{
\mathrm{ech}^{(2,\mathcal{A}^{(2)})}_{s}(\mathfrak{p}^{(2)})\colon
[M]_{s}
\ar@{=>}[r]^-{\text{\Small{($\mathfrak{p}^{(2)}$, $\mathrm{id}^{\mathrm{T}_{\Sigma^{\boldsymbol{\mathcal{A}}}}(X)_{s}}$)}}}
&
[N]_{s}
}
.
$$

Let us note that $\mathrm{ech}^{(2,\mathcal{A}^{(2)})}_{s}(\mathfrak{p}^{(2)})$ is neither a $(2,[1])$-identity second-order path, nor a second-order path of length strictly greater than one nor an echelonless second-order path. Therefore, for every sort $t\in S$ and every second-order path $\mathfrak{Q}^{(2)}\in\mathrm{Pth}_{\boldsymbol{\mathcal{A}}^{(2)},t}$, we have that $(\mathfrak{Q}^{(2)},t)\nprec_{\mathbf{Pth}_{\boldsymbol{\mathcal{A}}^{(2)}}} (\mathrm{ech}^{(2,\mathcal{A}^{(2)})}_{s}(\mathfrak{p}^{(2)}),s)$. Consequently, no element of $\coprod \mathrm{Pth}_{\boldsymbol{\mathcal{A}}^{(2)}}$ different from $(\mathrm{ech}^{(2,\mathcal{A}^{(2)})}_{s}(\mathfrak{p}^{(2)}),s)$ exists which $<_{\mathbf{Pth}_{\boldsymbol{\mathcal{A}}^{(2)}}}$-precedes $(\mathrm{ech}^{(2,\mathcal{A}^{(2)})}_{s}(\mathfrak{p}^{(2)}),s)$. From this it follows that $(\mathrm{ech}^{(2,\mathcal{A}^{(2)})}_{s}(\mathfrak{p}^{(2)}),s)$ is minimal in $(\coprod \mathrm{Pth}_{\boldsymbol{\mathcal{A}}^{(2)}}, \leq_{\mathbf{Pth}_{\boldsymbol{\mathcal{A}}^{(2)}}})$.

Conversely, we prove that every minimal element of $(\coprod \mathrm{Pth}_{\boldsymbol{\mathcal{A}}^{(2)}}, \leq_{\mathbf{Pth}_{\boldsymbol{\mathcal{A}}^{(2)}}})$ is a labelled $(2,[1])$-identity second-order path on a minimal path term class  or a labelled    second-order echelon.

Let $s$ be a sort in $S$ and $\mathfrak{P}^{(2)}$ a second-order path in $\mathrm{Pth}_{\boldsymbol{\mathcal{A}}^{(2)},s}$. Assume that $(\mathfrak{P}^{(2)},s)$ is a minimal element of $(\coprod \mathrm{Pth}_{\boldsymbol{\mathcal{A}}^{(2)}}, \leq_{\mathbf{Pth}_{\boldsymbol{\mathcal{A}}^{(2)}}})$. Then we want to show that $(\mathfrak{P}^{(2)},s)$ is in  $\coprod\mathrm{ip}^{(2,[1])\sharp}[
\mathrm{Min}(\coprod [\mathrm{PT}_{\boldsymbol{\mathcal{A}}}], \leq_{[\mathbf{PT}_{\boldsymbol{\mathcal{A}}}]})]$ or in the set $\coprod\mathrm{Ech}^{(2,\mathcal{A}^{(2)})}[\mathcal{A}^{(2)}]$.

We consider three cases. Either (1) $\mathfrak{P}^{(2)}$ is a $(2,[1])$-identity second-order path or (2) $\mathfrak{P}^{(2)}$ is a second-order path of length at least one containing echelons or (3) $\mathfrak{P}^{(2)}$ is an echelonless second-order path.

If (1) $\mathfrak{P}^{(2)}$ is a $(2,[1])$-identity second-order path of the form $\mathfrak{P}^{(2)}=\mathrm{ip}^{(2,[1])\sharp}_{s}([P]_{s})$, for a unique path term class $[P]_{s}$ in $[\mathrm{PT}_{\boldsymbol{\mathcal{A}}}]_{s}$, then either (1.1)  $([P]_{s},s)$ is a minimal pair or (1.2) $([P]_{s},s)$ is not a minimal pair in $(\coprod [\mathrm{PT}_{\boldsymbol{\mathcal{A}}}], \leq_{[\mathbf{PT}_{\boldsymbol{\mathcal{A}}}]})$.

If (1.1) $([P]_{s},s)$ is minimal in  $(\coprod [\mathrm{PT}_{\boldsymbol{\mathcal{A}}}], \leq_{[\mathbf{PT}_{\boldsymbol{\mathcal{A}}}]})$, then the statement follows because $(\mathfrak{P}^{(2)},s)$ is a labelled pair in $$\textstyle \coprod\mathrm{ip}^{(2,[1])\sharp}\left[
\mathrm{Min}\left(\coprod [\mathrm{PT}_{\boldsymbol{\mathcal{A}}}], \leq_{[\mathbf{PT}_{\boldsymbol{\mathcal{A}}}]}\right)\right].$$

If (1.2) $([P]_{s},s)$ is not minimal in  $(\coprod [\mathrm{PT}_{\boldsymbol{\mathcal{A}}}], \leq_{[\mathbf{PT}_{\boldsymbol{\mathcal{A}}}]})$, then one can find a sort $t\in S$ and a path class $[Q]_{t}$ in $[\mathrm{PT}_{\boldsymbol{\mathcal{A}}}]_{t}$ for which
$$
\left([Q]_{t}, t
\right)
<_{[\mathbf{PT}_{\boldsymbol{\mathcal{A}}}]}
\left([P]_{s}, s
\right).
$$

Consider the second-order path $\mathfrak{Q}^{(2)}$ in $\mathrm{Pth}_{\boldsymbol{\mathcal{A}}^{(2)},t}$ given by $\mathfrak{Q}^{(2)}=\mathrm{ip}^{(2,[1])\sharp}_{t}([Q]_{t})$. Then, taking into account the above strict inequality and according to Definition~\ref{DDOrd}, we have that 
$$
\left(\mathfrak{Q}^{(2)},t
\right)
<_{\mathbf{Pth}_{\boldsymbol{\mathcal{A}}^{(2)}}}
\left(\mathfrak{P}^{(2)},s
\right),
$$
contradicting that $(\mathfrak{P}^{(2)},s)$ is minimal.  This excludes case (1.2).

Now consider the case (2) in which  $\mathfrak{P}^{(2)}$ is a second-order path of length at least one containing   second-order echelons. Then either (2.1) the length of $\mathfrak{P}^{(2)}$ is one or (2.2) the length $\mathfrak{P}^{(2)}$ is strictly greater than one.

If (2.1), then $\mathfrak{P}^{(2)}$ is a   second-order echelon and the statement holds because $\mathfrak{P}^{(2)}$ is an element of $\mathrm{Ech}^{(2,\mathcal{A}^{(2)})}_{s}[\mathcal{A}^{(2)}_{s}]$.

If (2.2), then $\mathfrak{P}^{(2)}$ is a second-order path of length strictly greater than one containing at least one   second-order echelon. Then there exists an index $i\in\bb{\mathfrak{P}^{(2)}}$ such that the one-step subpath $\mathfrak{P}^{(2),i,i}$ is the first   second-order echelon appearing in $\mathfrak{P}^{(2)}$.  Regarding the nature of $i$, we have the following possibilities.

If  $i=0$ then, according to Definition~\ref{DDOrd}, we have that $(\mathfrak{P}^{(2),0,0},s)$ $\prec_{\mathbf{Pth}_{\boldsymbol{\mathcal{A}}^{(2)}}}$-precedes $
(\mathfrak{P}^{(2)},s)$. Let us note that $\mathfrak{P}^{(2),0,0}$ cannot be equal to $\mathfrak{P}^{(2)}$ because it has length $1$. This contradicts that $(\mathfrak{P}^{(2)},s)$ is minimal.

If  $i>0$, then, according to Definition~\ref{DDOrd}, we have that $(\mathfrak{P}^{(2),0,i-1},s)$ $\prec_{\mathbf{Pth}_{\boldsymbol{\mathcal{A}}^{(2)}}}$-precedes $
(\mathfrak{P}^{(2)},s)$. Let us note that $\mathfrak{P}^{(2),0,i-1}$ cannot be equal to $\mathfrak{P}^{(2)}$ because it has smaller length. This contradicts that $(\mathfrak{P}^{(2)},s)$ is minimal.  

In any case, this excludes second-order paths in case (2.2) of being minimal.

Let us consider, finally, the case (3) of $\mathfrak{P}^{(2)}$ being an echelonless second-order path. Then either (3.1) $\mathfrak{P}^{(2)}$ is not head-constant, or (3.2) $\mathfrak{P}^{(2)}$ is head-constant but not coherent or (3.3) $\mathfrak{P}^{(2)}$ is  coherent and head-constant.

If (3.1) $\mathfrak{P}^{(2)}$ is an echelonless second-order path that is not head-constant, then let $i\in \bb{\mathfrak{P}^{(2)}}$ be the maximum index for which $\mathfrak{P}^{(2),0,i}$ is head-constant then, according to Definition~\ref{DDOrd}, $(\mathfrak{P}^{(2),0,i},s)$ $\prec_{\mathbf{Pth}_{\boldsymbol{\mathcal{A}}^{(2)}}}$-precedes $
(\mathfrak{P}^{(2)},s)$. Let us note that $\mathfrak{P}^{(2),0,i}$ cannot be equal to $\mathfrak{P}^{(2)}$ because $\mathfrak{P}^{(2),0,i}$ is head-constant whilst $\mathfrak{P}^{(2)}$ is not head-constant.

If (3.2) $\mathfrak{P}^{(2)}$ is a head-constant echelonless second-order path that is not coherent, then let $i\in \bb{\mathfrak{P}^{(2)}}$ be the maximum index for which $\mathfrak{P}^{(2),0,i}$ is coherent then, according to Definition~\ref{DDOrd}, $(\mathfrak{P}^{(2),0,i},s)$ $\prec_{\mathbf{Pth}_{\boldsymbol{\mathcal{A}}^{(2)}}}$-precedes $
(\mathfrak{P}^{(2)},s)$. Let us note that $\mathfrak{P}^{(2),0,i}$ cannot be equal to $\mathfrak{P}^{(2)}$ because $\mathfrak{P}^{(2),0,i}$ is coherent whilst $\mathfrak{P}^{(2)}$ is not coherent. 

If (3.3) $\mathfrak{P}^{(2)}$ is a head-constant coherent echelonless second-order path, then this second-order path is suitable for the application of the second-order extraction procedure. Let $\mathfrak{Q}^{(2)}$ be a second-order path of sort $t\in S$ extracted from $\mathfrak{P}^{(2)}$ according to Lemma~\ref{LDPthExtract} then, according to Definition~\ref{DDOrd}, $(\mathfrak{Q}^{(2)},t)$ $\prec_{\mathbf{Pth}_{\boldsymbol{\mathcal{A}}^{(2)}}}$-precedes $
(\mathfrak{P}^{(2)},s)$. It follows from Lemma~\ref{LDPthExtract} that $\mathfrak{Q}^{(2)}$ cannot be equal to $\mathfrak{P}^{(2)}$. 

In any case, we have found a labelled second-order path different from $(\mathfrak{P}^{(2)},s)$ that $\leq_{\mathbf{Pth}_{\boldsymbol{\mathcal{A}}^{(2)}}}$-precedes it. This excludes second-order paths in case (3) of being minimal.
\end{proof}

In the following proposition we prove that the coproduct of the many-sorted set of $(2,[1])$-identity second-order paths is a lower set of $(\coprod\mathrm{Pth}_{\boldsymbol{\mathcal{A}}^{(2)}}, \leq_{\mathbf{Pth}_{\boldsymbol{\mathcal{A}}^{(2)}}})$. Moreover this subset reflects the preorder arising from $[\mathrm{PT}_{\boldsymbol{\mathcal{A}}}]$.

\begin{restatable}{proposition}{PDULower}
\label{PDULower} The subset $\coprod\mathrm{ip}^{(2,[1])\sharp}[[\mathrm{PT}_{\boldsymbol{\mathcal{A}}}]]$ of $\coprod\mathrm{Pth}_{\boldsymbol{\mathcal{A}}^{(2)}}$, consisting of all ordered pairs of labeled $(2,[1])$-identity second-order paths, 
is a lower set of the partially ordered set $(\coprod\mathrm{Pth}_{\boldsymbol{\mathcal{A}}^{(2)}}, \leq_{\mathbf{Pth}_{\boldsymbol{\mathcal{A}}^{(2)}}})$. That is, for every pair of sorts $s,t\in S$ and second-order paths $\mathfrak{P}^{(2)}\in\mathrm{Pth}_{\boldsymbol{\mathcal{A}}^{(2)},s}$ and $\mathfrak{Q}^{(2)}\in\mathrm{Pth}_{\boldsymbol{\mathcal{A}}^{(2)},t}$ if
\begin{enumerate}
\item $(\mathfrak{Q}^{(2)},t)\leq_{\mathbf{Pth}_{\boldsymbol{\mathcal{A}}^{(2)}}}(\mathfrak{P}^{(2)},s)$ and
\item $\mathfrak{P}^{(2)}$ is a $(2,[1])$-identity second-order path,
\end{enumerate}
then $\mathfrak{Q}^{(2)}$ is also a $(2,[1])$-identity second-order path. Moreover, if $\mathfrak{P}^{(2)}$ and $\mathfrak{Q}^{(2)}$ have the form
\begin{align*}
\mathfrak{P}^{(2)}&=
\mathrm{ip}^{(2,[1])\sharp}_{s}\left(
[P]_{s}
\right),
&
\mathfrak{Q}^{(2)}&=
\mathrm{ip}^{(2,[1])\sharp}_{t}\left(
[Q]_{t}
\right),
\end{align*}
for some path term classes $[P]_{s}\in[\mathrm{PT}_{\boldsymbol{\mathcal{A}}}]_{s}$ and $[Q]_{t}\in[\mathrm{PT}_{\boldsymbol{\mathcal{A}}}]_{t}$, then the following inequality holds
$$
\left([Q]_{t},t\right)
\leq_{[\mathbf{PT}_{\boldsymbol{\mathcal{A}}}]}
\left([P]_{s},s\right).
$$
\end{restatable}

\begin{proof}
Let $(\mathfrak{Q}^{(2)},t)$, $(\mathfrak{P}^{(2)},s)$ be two pairs in $\coprod\mathrm{Pth}_{\boldsymbol{\mathcal{A}}^{(2)}}$ with $(\mathfrak{Q}^{(2)},t)\leq_{\mathbf{Pth}_{\boldsymbol{\mathcal{A}}^{(2)}}}(\mathfrak{P}^{(2)},s)$. Let us assume that $\mathfrak{P}^{(2)}$ is a $(2,[1])$-identity second-order path. 
Since $(\mathfrak{Q}^{(2)},t)\leq_{\mathbf{Pth}_{\boldsymbol{\mathcal{A}}^{(2)}}}(\mathfrak{P}^{(2)},s)$, we have, by Remark~\ref{RDOrd}, that either $s=t$ and $\mathfrak{Q}^{(2)}=\mathfrak{P}^{(2)}$ or there exists a natural number $m\in\mathbb{N}-\{0\}$, a word $\mathbf{w}\in S^{\star}$ of length $\bb{\mathbf{w}}=m+1$, and a family of second-order paths $(\mathfrak{R}^{(2)}_{k})_{k\in\bb{\mathbf{w}}}$ in $\mathrm{Pth}_{\boldsymbol{\mathcal{A}}^{(2)},\mathbf{w}}$ such that $w_{0}=t$, $\mathfrak{R}^{(2)}_{0}=\mathfrak{Q}^{(2)}$, $w_{m}=s$, $\mathfrak{R}^{(2)}_{m}=\mathfrak{P}^{(2)}$ and, for every $k\in m$, 
$
(\mathfrak{R}^{(2)}_{k}, w_{k})
\prec_{\mathbf{Pth}_{\boldsymbol{\mathcal{A}}^{(2)}}}
(\mathfrak{R}^{(2)}_{k+1}, w_{k+1}).
$

If $s=t$ and $\mathfrak{Q}^{(2)}=\mathfrak{P}^{(2)}$ the statement trivially holds. Therefore, it suffices to prove the proposition for a non-trivial sequence instantiating that 
$$\left(\mathfrak{Q}^{(2)},t
\right)\leq_{\mathbf{Pth}_{\boldsymbol{\mathcal{A}}^{(2)}}}
\left(\mathfrak{P}^{(2)},s
\right).$$

The proof is by induction on $m\in\mathbb{N}-\{0\}$.

\textsf{Base step of the induction.}

For $m=1$ we have that $(\mathfrak{Q}^{(2)},t)\prec_{\mathbf{Pth}_{\boldsymbol{\mathcal{A}}^{(2)}}}(\mathfrak{P}^{(2)},s)$. Then, according to Definition~\ref{DDOrd}, since $\mathfrak{P}^{(2)}$ is a $(2,[1])$-identity second-order path, we have that the unique possibility for $\mathfrak{Q}^{(2)}$ is that it is also a $(2,[1])$-identity second-order path. 

Moreover, if the second-order paths $\mathfrak{Q}^{(2)}$ and $\mathfrak{P}^{(2)}$ have the form
\begin{align*}
\mathfrak{P}^{(2)}&=
\mathrm{ip}^{(2,[1])\sharp}_{s}\left(
[P]_{s}
\right),
&
\mathfrak{Q}^{(2)}&=
\mathrm{ip}^{(2,[1])\sharp}_{t}
\left(
[Q]_{t}
\right),
\end{align*}
for some path term classes $[P]_{s}\in[\mathrm{PT}_{\boldsymbol{\mathcal{A}}}]_{s}$ and $[Q]_{t}\in[\mathrm{PT}_{\boldsymbol{\mathcal{A}}}]_{t}$, then the following inequality also holds
$$
\left([Q]_{t},t
\right)
\leq_{[\mathbf{PT}_{\boldsymbol{\mathcal{A}}}]}
\left([P]_{s},s
\right).
$$

\textsf{Inductive step of the induction.}

Let us suppose that the statement holds for sequences of length up to $m\in\mathbb{N}-\{0\}$. That is, if $\mathbf{w}$ is a word in $S^{\star}$ of length $\bb{\mathbf{w}}=m+1$ and $(\mathfrak{R}^{(2)}_{k})_{k\in\bb{\mathbf{w}}}$ is a family of second-order paths in $\mathrm{Pth}_{\boldsymbol{\mathcal{A}}^{(2)},\mathbf{w}}$ such that $w_{0}=t$, $\mathfrak{R}^{(2)}_{0}=\mathfrak{Q}^{(2)}$, $w_{m}=s$, $\mathfrak{R}^{(2)}_{m}=\mathfrak{P}^{(2)}$ and, for every $k\in m$, $(\mathfrak{R}^{(2)}_{k},w_{k})\prec_{\mathbf{Pth}_{\boldsymbol{\mathcal{A}}^{(2)}}} (\mathfrak{R}^{(2)}_{k+1}, w_{k+1})$ then, if $\mathfrak{P}^{(2)}$ is a $(2,[1])$-identity second-order path, it follows that $\mathfrak{Q}^{(2)}$ is also a $(2,[1])$-identity second-order path. 

Moreover, if the second-order paths $\mathfrak{Q}^{(2)}$ and $\mathfrak{P}^{(2)}$ have the form
\begin{align*}
\mathfrak{P}^{(2)}&=
\mathrm{ip}^{(2,[1])\sharp}_{s}\left(
[P]_{s}
\right),
&
\mathfrak{Q}^{(2)}&=
\mathrm{ip}^{(2,[1])\sharp}_{t}\left(
[Q]_{t}
\right),
\end{align*}
for some path term classes $[P]_{s}\in[\mathrm{PT}_{\boldsymbol{\mathcal{A}}}]_{s}$ and $[Q]_{t}\in[\mathrm{PT}_{\boldsymbol{\mathcal{A}}}]_{t}$, then the following inequality also holds
$$
\left([Q]_{t},t\right)
\leq_{[\mathbf{PT}_{\boldsymbol{\mathcal{A}}}]}
\left([P]_{s},s\right).
$$

Now, let $\mathbf{w}$ be a word in $S^{\star}$ of length $\bb{\mathbf{w}}=m+2$ and $(\mathfrak{R}^{(2)}_{k})_{k\in\bb{\mathbf{w}}}$ a family of second-order paths in $\mathrm{Pth}_{\boldsymbol{\mathcal{A}}^{(2)},\mathbf{w}}$ such that $w_{0}=t$, $\mathfrak{R}^{(2)}_{0}=\mathfrak{Q}^{(2)}$, $w_{m+1}=s$, $\mathfrak{R}^{(2)}_{m+1}=\mathfrak{P}^{(2)}$ and, for every $k\in m+1$, $(\mathfrak{R}^{(2)}_{k},w_{k})\prec_{\mathbf{Pth}_{\boldsymbol{\mathcal{A}}^{(2)}}} (\mathfrak{R}^{(2)}_{k+1}, w_{k+1})$. 

Consider the final subsequence $(\mathfrak{R}^{(2)}_{k})_{k\in [1,m+1]}$. This is a sequence of length $m$ instantiating that $(\mathfrak{R}^{(2)}_{1},w_{1})\leq_{\mathbf{Pth}_{\boldsymbol{\mathcal{A}}^{(2)}}} (\mathfrak{P}^{(2)},s)$. Now, since $\mathfrak{P}^{(2)}$ is a $(2,[1])$-identity second-order path, it follows, by induction, that $\mathfrak{R}^{(2)}_{1}$ is also a $(2,[1])$-identity second-order path. 

Moreover, if the second-order paths $\mathfrak{R}^{(2)}_{1}$ and $\mathfrak{P}^{(2)}$ have the form
\begin{align*}
\mathfrak{P}^{(2)}&=
\mathrm{ip}^{(2,[1])\sharp}_{s}\left(
[P]_{s}
\right),
&
\mathfrak{R}^{(2)}_{1}&=
\mathrm{ip}^{(2,[1])\sharp}_{w_{1}}\left(
[R_{1}]_{w_{1}}
\right),
\end{align*}
for some path term classes $[P]_{s}\in[\mathrm{PT}_{\boldsymbol{\mathcal{A}}}]_{s}$ and $[R_{1}]_{w_{1}}\in[\mathrm{PT}_{\boldsymbol{\mathcal{A}}}]_{w_{1}}$, then the following inequality also holds
$$
\left([R_{1}]_{w_{1}},w_{1}
\right)
\leq_{[\mathbf{PT}_{\boldsymbol{\mathcal{A}}}]}
\left([P]_{s},s
\right).
$$

Now, consider the initial subsequence $(\mathfrak{R}^{(2)}_{k})_{k\in 2}$. This is a sequence of length $1$ instantiating that $(\mathfrak{Q}^{(2)},t)\prec_{\mathbf{Pth}_{\boldsymbol{\mathcal{A}}^{(2)}}} (\mathfrak{R}^{(2)},w_{1})$. Now, since $\mathfrak{R}^{(2)}$ is a $(2,[1])$-identity second-order path, it follows, by the base step, that $\mathfrak{Q}^{(2)}$ is also a $(2,[1])$-identity second-order path. 

Moreover, if the second-order paths $\mathfrak{R}^{(2)}_{1}$ and $\mathfrak{Q}^{(2)}$ have the form
\begin{align*}
\mathfrak{Q}^{(2)}&=
\mathrm{ip}^{(2,[1])\sharp}_{t}\left(
[Q]_{t}
\right),
&
\mathfrak{R}^{(2)}_{1}&=
\mathrm{ip}^{(2,[1])\sharp}_{w_{1}}\left(
[R_{1}]_{w_{1}}
\right),
\end{align*}
for some path term classes $[Q]_{t}\in[\mathrm{PT}_{\boldsymbol{\mathcal{A}}}]_{t}$ and $[R_{1}]_{w_{1}}\in[\mathrm{PT}_{\boldsymbol{\mathcal{A}}}]_{w_{1}}$, then the following inequality also holds
$$
\left([Q]_{t},t
\right)
\leq_{[\mathbf{PT}_{\boldsymbol{\mathcal{A}}}]}
\left([R_{1}]_{w_{1}},w_{1}
\right).
$$

All in all, we can affirm that $\mathfrak{Q}^{(2)}$ is a $(2,[1])$-identity second-order path and, taking into account Proposition~\ref{PPTQOrdArt}, since $(\coprod[\mathrm{PT}_{\boldsymbol{\mathcal{A}}}], \leq_{[\mathbf{PT}_{\boldsymbol{\mathcal{A}}}]})$ is a partially ordered set, we have, by transitivity, that the following inequality holds
$$
\left([Q]_{t},t
\right)
\leq_{[\mathbf{PT}_{\boldsymbol{\mathcal{A}}}]}
\left([P]_{s},s
\right).
$$

This finishes the proof.
\end{proof}

From the last proposition we obtain that every strictly decreasing sequence instantiating an inequality where the greatest element is a $(2,[1])$-identity second-order path is entirely composed of $(2,[1])$-identity second-order paths. 

\begin{corollary}\label{CDULower}
Let $m$ be a natural number in $\mathbb{N}-\{0\}$, $\mathbf{w}$ a word in $S^{\star}$ such that $\bb{\mathbf{w}}=m+1$ and $(\mathfrak{R}^{(2)}_{k})_{k\in \bb{\mathbf{w}}}$ a family of second-order paths in $\mathrm{Pth}_{\boldsymbol{\mathcal{A}}^{(2)},\mathbf{w}}$ such that
\begin{enumerate}
\item for every $k\in m$, $(\mathfrak{R}^{(2)}_{k}, w_{k})\prec_{\mathbf{Pth}_{\boldsymbol{\mathcal{A}}^{(2)}}} (\mathfrak{R}^{(2)}_{k+1}, w_{k+1})$ and 
\item $\mathfrak{R}^{(2)}_{m}$ is a $(2,[1])$-identity second-order path.
\end{enumerate}
Then, for every $k\in m$, $\mathfrak{R}^{(2)}_{k}$ is a $(2,[1])$-identity second-order path. Moreover, for every $k\in m+1$, if $\mathfrak{R}^{(2)}_{k}$ has the form
\begin{align*}
\mathfrak{R}^{(2)}_{k}&=
\mathrm{ip}^{(2,[1])\sharp}_{w_{k}}\left(
[R_{k}]_{w_{k}}
\right),
\end{align*}
for some path term class $[R_{k}]_{w_{k}}\in[\mathrm{PT}_{\boldsymbol{\mathcal{A}}}]_{w_{k}}$, then, for every $k\in m$, the following inequality  holds
$$
\left([R_{k}]_{w_{k}},w_{k}\right)
\leq_{[\mathbf{PT}_{\boldsymbol{\mathcal{A}}}]}
\left([R_{k+1}]_{w_{k+1}},w_{k+1}
\right).
$$
\end{corollary}

We next prove some technical results concerning the nature of pairs in the preorder $\leq_{\mathbf{Pth}_{\boldsymbol{\mathcal{A}}^{(2)}}}$. In this respect, we will show that, for every pair of sorts $s$ and $t$ in $S$ and every pair of second-order paths $\mathfrak{P}^{(2)}\in\mathrm{Pth}_{\boldsymbol{\mathcal{A}}^{(2)},s}$, $\mathfrak{Q}^{(2)}\in\mathrm{Pth}_{\boldsymbol{\mathcal{A}}^{(2)},t}$, if $(\mathfrak{Q}^{(2)},t)<_{\mathbf{Pth}_{\boldsymbol{\mathcal{A}}^{(2)}}}(\mathfrak{P}^{(2)},s)$, then $\mathfrak{P}^{(2)}$ is strictly more complex than $\mathfrak{Q}^{(2)}$ with respect to either length or height of the associated first-order translation.

\begin{lemma}\label{LDOrdI} Let $m$ be a natural number in $\mathbb{N}-\{0\}$, $\mathbf{w}$ a word in $S^{\star}$ such that $\bb{\mathbf{w}}=m+1$ and $(\mathfrak{R}^{(2)}_{k})_{k\in \bb{\mathbf{w}}}$ a family of second-order paths in $\mathrm{Pth}_{\boldsymbol{\mathcal{A}}^{(2)},\mathbf{w}}$ such that, for every $k\in m$, $(\mathfrak{R}^{(2)}_{k}, w_{k})\prec_{\mathbf{Pth}_{\boldsymbol{\mathcal{A}}^{(2)}}} (\mathfrak{R}^{(2)}_{k+1}, w_{k+1})$. Then, for every $k\in [1,m]$, we have that
\begin{enumerate}
\item $\mathfrak{R}^{(2)}_{k}$ is a $(2,[1])$-identity second-order path on a non-minimal path term class; or
\item $\mathfrak{R}^{(2)}_{k}$ is a second-order path of length strictly greater than one containing at least one   second-order echelon; or
\item $\mathfrak{R}^{(2)}_{k}$ is an echelonless second-order path.
\end{enumerate}
\end{lemma}

\begin{proof}
Let us recall (see Remark~\ref{RDIrrefl}) that $\prec_{\mathbf{Pth}_{\boldsymbol{\mathcal{A}}^{(2)}}}$ and $\Delta_{\coprod\mathrm{Pth}_{\boldsymbol{\mathcal{A}}^{(2)}}}$ are disjoint sets. Moreover, by hypothesis, for every $k\in m$, $(\mathfrak{R}^{(2)}_{k}, w_{k})\prec_{\mathbf{Pth}_{\boldsymbol{\mathcal{A}}^{(2)}}} (\mathfrak{R}^{(2)}_{k+1}, w_{k+1})$. From this it follows that, for every $k\in [1,m]$, the element $(\mathfrak{R}^{(2)}_{k}, w_{k})$ is not a minimal element of $(\coprod \mathrm{Pth}_{\boldsymbol{\mathcal{A}}^{(2)}}, \leq_{\mathbf{Pth}_{\boldsymbol{\mathcal{A}}^{(2)}}})$.
\end{proof}

In the following lemma we show that, for every pair of sorts $s$ and $t$ in $S$ and every pair of second-order paths $\mathfrak{P}^{(2)}\in\mathrm{Pth}_{\boldsymbol{\mathcal{A}}^{(2)},s}$ and $\mathfrak{Q}^{(2)}\in\mathrm{Pth}_{\boldsymbol{\mathcal{A}}^{(2)},t}$, if $(\mathfrak{Q}^{(2)},t)<_{\mathbf{Pth}_{\boldsymbol{\mathcal{A}}^{(2)}}}(\mathfrak{P}^{(2)},s)$ and   there exists an strictly decreasing sequence instantiating this inequality composed of  second-order paths of length strictly greater than one containing at least one   second-order echelon, then $\bb{\mathfrak{Q}^{(2)}}<\bb{\mathfrak{P}^{(2)}}$.

\begin{lemma}\label{LDOrdII} 
Let $m$ be a natural number in $\mathbb{N}-\{0\}$, $\mathbf{w}$ a word in $S^{\star}$ such that $\bb{\mathbf{w}}=m+1$ and  $(\mathfrak{R}^{(2)}_{k})_{k\in \bb{\mathbf{w}}}$ a family of second-order paths in $\mathrm{Pth}_{\boldsymbol{\mathcal{A}}^{(2)},\mathbf{w}}$ such that
\begin{enumerate}
\item for every $k\in m$, $(\mathfrak{R}^{(2)}_{k}, w_{k})\prec_{\mathbf{Pth}_{\boldsymbol{\mathcal{A}}^{(2)}}} (\mathfrak{R}^{(2)}_{k+1}, w_{k+1})$ and, 
\item for every $k\in [1,m]$, $\mathfrak{R}^{(2)}_{k}$ is a second-order path of length strictly greater than one containing at least one   second-order echelon.
\end{enumerate}
Then
$\bb{\mathfrak{R}^{(2)}_{0}}<\bb{\mathfrak{R}^{(2)}_{m}}$.
\end{lemma}
\begin{proof}
We prove it by induction on $m\in\mathbb{N}-\{0\}$.

\textsf{Base step of the induction.}

For $m=1$, we have that $(\mathfrak{R}^{(2)}_{0},w_{0})\prec_{\mathbf{Pth}_{\boldsymbol{\mathcal{A}}^{(2)}}}(\mathfrak{R}^{(2)}_{1},w_{1})$, where $\mathfrak{R}^{(2)}_{1}$ is a second-order path of length strictly greater than one containing at least one   second-order echelon.

Let $i\in\bb{\mathfrak{R}^{(2)}_{1}}$ be the first index for which $\mathfrak{R}^{(2),i,i}_{1}$ is a   second-order echelon. Then, according to Definition~\ref{DDOrd}, we have that $w_{0}=w_{1}$ and, taking into account the different possibilities for $i$, we have that $\mathfrak{R}^{(2)}_{0}$ has to be equal to one of the following subpaths of $\mathfrak{R}^{(2)}_{1}$
\begin{align*}
\mathfrak{R}^{(2),0,0}_{1},&&
\mathfrak{R}^{(2),1,\bb{\mathfrak{R}^{(2)}_{1}}-1}_{1},&&
\mathfrak{R}^{(2),0,i-1}_{1}, \qquad \mbox{or}&&
\mathfrak{R}^{(2),i,\bb{\mathfrak{R}^{(2)}_{1}}-1}_{1}.
\end{align*}
Note that, in every case, we can guarantee that  $\bb{\mathfrak{R}^{(2)}_{0}}<\bb{\mathfrak{R}^{(2)}_{1}}$.

\textsf{Inductive step of the induction.}

Let us suppose that the statement holds for sequences of length up to $m\in\mathbb{N}-\{0\}$. That is, if $\mathbf{w}$ is a word in $S^{\star}$ of length $\bb{\mathbf{w}}=m+1$ and $(\mathfrak{R}^{(2)}_{k})_{k\in\bb{\mathbf{w}}}$ is a family of second-order paths in $\mathrm{Pth}_{\boldsymbol{\mathcal{A}}^{(2)},\mathbf{w}}$ satisfying that 
\begin{enumerate}
\item for every $k\in m$, $(\mathfrak{R}^{(2)}_{k}, w_{k})\prec_{\mathbf{Pth}_{\boldsymbol{\mathcal{A}}^{(2)}}} (\mathfrak{R}^{(2)}_{k+1}, w_{k+1})$; and
\item for every $k\in [1,m]$, $\mathfrak{R}^{(2)}_{k}$ is a second-order path of length strictly greater than one containing at least one   second-order echelon,
\end{enumerate}
then
$\bb{\mathfrak{R}^{(2)}_{0}}<\bb{\mathfrak{R}^{(2)}_{m}}$.

Now, let $\mathbf{w}$ be a word in $S^{\star}$ of length $\bb{\mathbf{w}}=m+2$ and let $(\mathfrak{R}^{(2)}_{k})_{k\in\bb{\mathbf{w}}}$ be a family of second-order paths in $\mathrm{Pth}_{\boldsymbol{\mathcal{A}}^{(2)},\mathbf{w}}$ satisfying that 
\begin{enumerate}
\item for every $k\in m$, $(\mathfrak{R}^{(2)}_{k}, w_{k})\prec_{\mathbf{Pth}_{\boldsymbol{\mathcal{A}}^{(2)}}} (\mathfrak{R}^{(2)}_{k+1}, w_{k+1})$; and
\item for every $k\in [1,m]$, $\mathfrak{R}^{(2)}_{k}$ is a second-order path of length strictly greater than one containing at least one second-order echelon.
\end{enumerate}

Consider the final subsequence $(\mathfrak{R}^{(2)}_{k})_{k\in [1,m+1]}$. This is a sequence of length $m$ satisfying the conditions of the inductive step. Hence, $\bb{\mathfrak{R}^{(2)}_{1}}<\bb{\mathfrak{R}^{(2)}_{m+1}}$.

Now, consider the initial subsequence $(\mathfrak{R}^{(2)}_{k})_{k\in 2}$. This is a sequence of length $1$ instantiating that $(\mathfrak{R}^{(2)}_{0},w_{0})\prec_{\mathbf{Pth}_{\boldsymbol{\mathcal{A}}^{(2)}}}(\mathfrak{R}^{(2)}_{1},w_{1})$ and, by assumption, $\mathfrak{R}^{(2)}_{1}$ is a second-order path of length strictly greater than one containing at least one   second-order echelon. By the base step, we have that $\bb{\mathfrak{R}^{(2)}_{0}}<\bb{\mathfrak{R}^{(2)}_{1}}$.

All in all, we can affirm that $\bb{\mathfrak{R}^{(2)}_{0}}<\bb{\mathfrak{R}^{(2)}_{m+1}}$.

This finishes the proof.
\end{proof}

In the following lemma we show that, for every pair of sorts $s$ and $t$ in $S$ and every pair of second-order paths $\mathfrak{P}^{(2)}\in\mathrm{Pth}_{\boldsymbol{\mathcal{A}}^{(2)},s}$ and $\mathfrak{Q}^{(2)}\in\mathrm{Pth}_{\boldsymbol{\mathcal{A}}^{(2)},t}$, if $(\mathfrak{Q}^{(2)},t)<_{\mathbf{Pth}_{\boldsymbol{\mathcal{A}}^{(2)}}}(\mathfrak{P}^{(2)},s)$ and   there exists an strictly decreasing sequence instantiating this inequality composed of  echelonless second-order paths that are not head-constant, then $\bb{\mathfrak{Q}^{(2)}}<\bb{\mathfrak{P}^{(2)}}$.

\begin{lemma}\label{LDOrdIII} 
Let $m$ be a natural number in $\mathbb{N}-\{0\}$, $\mathbf{w}$ a word in $S^{\star}$ such that $\bb{\mathbf{w}}=m+1$ and  $(\mathfrak{R}^{(2)}_{k})_{k\in \bb{\mathbf{w}}}$ a family of second-order paths in $\mathrm{Pth}_{\boldsymbol{\mathcal{A}}^{(2)},\mathbf{w}}$ such that
\begin{enumerate}
\item for every $k\in m$, $(\mathfrak{R}^{(2)}_{k}, w_{k})\prec_{\mathbf{Pth}_{\boldsymbol{\mathcal{A}}^{(2)}}} (\mathfrak{R}^{(2)}_{k+1}, w_{k+1})$ and, 
\item for every $k\in [1,m]$, $\mathfrak{R}^{(2)}_{k}$ is an echelonless second-order path  that is not head-constant.
\end{enumerate}
Then
$\bb{\mathfrak{R}^{(2)}_{0}}<\bb{\mathfrak{R}^{(2)}_{m}}$.
\end{lemma}
\begin{proof}
We prove it by induction on $m\in\mathbb{N}-\{0\}$.

\textsf{Base step of the induction.}

For $m=1$, we have that $(\mathfrak{R}^{(2)}_{0},w_{0})\prec_{\mathbf{Pth}_{\boldsymbol{\mathcal{A}}^{(2)}}}(\mathfrak{R}^{(2)}_{1},w_{1})$, where $\mathfrak{R}^{(2)}_{1}$ is an echelonless second-order path that is not head-constant.

Let $i\in\bb{\mathfrak{R}^{(2)}_{1}}$ be the maximum index for which $\mathfrak{R}^{(2),0,i}_{1}$ is a head-constant echelonless second-order path. Then, according to Definition~\ref{DDOrd}, we have that $w_{0}=w_{1}$ and $\mathfrak{R}^{(2)}_{0}$ has to be equal to one of the following subpaths of $\mathfrak{R}^{(2)}_{1}$
\begin{align*}
\mathfrak{R}^{(2),0,i}_{1},&&
\mathfrak{R}^{(2),i+1,\bb{\mathfrak{R}^{(2)}_{1}}-1}_{1}.
\end{align*}
Note that, in every case, we can guarantee that  $\bb{\mathfrak{R}^{(2)}_{0}}<\bb{\mathfrak{R}^{(2)}_{1}}$.

\textsf{Inductive step of the induction.}

Let us suppose that the statement holds for sequences of length up to $m\in\mathbb{N}-\{0\}$. That is, if $\mathbf{w}$ is a word in $S^{\star}$ of length $\bb{\mathbf{w}}=m+1$ and $(\mathfrak{R}^{(2)}_{k})_{k\in\bb{\mathbf{w}}}$ is a family of second-order paths in $\mathrm{Pth}_{\boldsymbol{\mathcal{A}}^{(2)},\mathbf{w}}$ satisfying that 
\begin{enumerate}
\item for every $k\in m$, $(\mathfrak{R}^{(2)}_{k}, w_{k})\prec_{\mathbf{Pth}_{\boldsymbol{\mathcal{A}}^{(2)}}} (\mathfrak{R}^{(2)}_{k+1}, w_{k+1})$; and
\item for every $k\in [1,m]$, $\mathfrak{R}^{(2)}_{k}$ is an echelonless second-order path that is not head-constant,
\end{enumerate}
then
$\bb{\mathfrak{R}^{(2)}_{0}}<\bb{\mathfrak{R}^{(2)}_{m}}$.

Now, let $\mathbf{w}$ be a word in $S^{\star}$ of length $\bb{\mathbf{w}}=m+2$ and let $(\mathfrak{R}^{(2)}_{k})_{k\in\bb{\mathbf{w}}}$ be a family of second-order paths in $\mathrm{Pth}_{\boldsymbol{\mathcal{A}}^{(2)},\mathbf{w}}$ satisfying that 
\begin{enumerate}
\item for every $k\in m$, $(\mathfrak{R}^{(2)}_{k}, w_{k})\prec_{\mathbf{Pth}_{\boldsymbol{\mathcal{A}}^{(2)}}} (\mathfrak{R}^{(2)}_{k+1}, w_{k+1})$; and
\item for every $k\in [1,m]$, $\mathfrak{R}^{(2)}_{k}$ is an echelonless second-order path  that is not head-constant.
\end{enumerate}

Consider the final subsequence $(\mathfrak{R}^{(2)}_{k})_{k\in [1,m+1]}$. This is a sequence of length $m$ satisfying the conditions of the inductive step. Hence, $\bb{\mathfrak{R}^{(2)}_{1}}<\bb{\mathfrak{R}^{(2)}_{m+1}}$.

Now, consider the initial subsequence $(\mathfrak{R}^{(2)}_{k})_{k\in 2}$. This is a sequence of length $1$ instantiating that $(\mathfrak{R}^{(2)}_{0},w_{0})\prec_{\mathbf{Pth}_{\boldsymbol{\mathcal{A}}^{(2)}}}(\mathfrak{R}^{(2)}_{1},w_{1})$ and, by assumption, $\mathfrak{R}^{(2)}_{1}$ is an echelonless second-order path that is not head-constant. By the base step, we have that $\bb{\mathfrak{R}^{(2)}_{0}}<\bb{\mathfrak{R}^{(2)}_{1}}$.

All in all, we can affirm that $\bb{\mathfrak{R}^{(2)}_{0}}<\bb{\mathfrak{R}^{(2)}_{m+1}}$.

This finishes the proof.
\end{proof}

In the following lemma we show that, for every pair of sorts $s$ and $t$ in $S$ and every pair of second-order paths $\mathfrak{P}^{(2)}\in\mathrm{Pth}_{\boldsymbol{\mathcal{A}}^{(2)},s}$ and $\mathfrak{Q}^{(2)}\in\mathrm{Pth}_{\boldsymbol{\mathcal{A}}^{(2)},t}$, if $(\mathfrak{Q}^{(2)},t)<_{\mathbf{Pth}_{\boldsymbol{\mathcal{A}}^{(2)}}}(\mathfrak{P}^{(2)},s)$ and there exists an strictly decreasing sequence instantiating this inequality composed of head-constant echelonless second-order paths that are not coherent, then $\bb{\mathfrak{Q}^{(2)}}<\bb{\mathfrak{P}^{(2)}}$.

\begin{lemma}\label{LDOrdIV} 
Let $m$ be a natural number in $\mathbb{N}-\{0\}$, $\mathbf{w}$ a word in $S^{\star}$ such that $\bb{\mathbf{w}}=m+1$, and  $(\mathfrak{R}^{(2)}_{k})_{k\in \bb{\mathbf{w}}}$ a family of second-order paths in $\mathrm{Pth}_{\boldsymbol{\mathcal{A}}^{(2)},\mathbf{w}}$ such that
\begin{enumerate}
\item for every $k\in m$, $(\mathfrak{R}^{(2)}_{k}, w_{k})\prec_{\mathbf{Pth}_{\boldsymbol{\mathcal{A}}^{(2)}}} (\mathfrak{R}^{(2)}_{k+1}, w_{k+1})$; and
\item for every $k\in [1,m]$, $\mathfrak{R}^{(2)}_{k}$ is a head-constant echelonless second-order path that is not coherent.
\end{enumerate}
Then
$\bb{\mathfrak{R}^{(2)}_{0}}<\bb{\mathfrak{R}^{(2)}_{m}}$.
\end{lemma}
\begin{proof}
We prove it by induction on $m\in\mathbb{N}-\{0\}$.

\textsf{Base step of the induction.}

For $m=1$, we have that $(\mathfrak{R}^{(2)}_{0},w_{0})\prec_{\mathbf{Pth}_{\boldsymbol{\mathcal{A}}^{(2)}}}(\mathfrak{R}^{(2)}_{1},w_{1})$, where $\mathfrak{R}^{(2)}_{1}$ is a head-constant echelonless second-order path that is not coherent.

Let $i\in\bb{\mathfrak{R}^{(2)}_{1}}$ be the maximum index for which $\mathfrak{R}^{(2),0,i}_{1}$ is a coherent head-constant echelonless second-order path. Then, according to Definition~\ref{DDOrd}, we have that $w_{0}=w_{1}$ and $\mathfrak{R}^{(2)}_{0}$ has to be equal to one of the following subpaths of $\mathfrak{R}^{(2)}_{1}$
\begin{align*}
\mathfrak{R}^{(2),0,i}_{1},&&
\mathfrak{R}^{(2),i+1,\bb{\mathfrak{R}^{(2)}_{1}}-1}_{1}.
\end{align*}
Note that, in every case, we can guarantee that  $\bb{\mathfrak{R}^{(2)}_{0}}<\bb{\mathfrak{R}^{(2)}_{1}}$.

\textsf{Inductive step of the induction.}

Let us suppose that the statement holds for sequences of length up to $m\in\mathbb{N}-\{0\}$. That is, if $\mathbf{w}$ is a word in $S^{\star}$ of length $\bb{\mathbf{w}}=m+1$ and $(\mathfrak{R}^{(2)}_{k})_{k\in\bb{\mathbf{w}}}$ is a family of second-order paths in $\mathrm{Pth}_{\boldsymbol{\mathcal{A}}^{(2)},\mathbf{w}}$ satisfying that 
\begin{enumerate}
\item for every $k\in m$, $(\mathfrak{R}^{(2)}_{k}, w_{k})\prec_{\mathbf{Pth}_{\boldsymbol{\mathcal{A}}^{(2)}}} (\mathfrak{R}^{(2)}_{k+1}, w_{k+1})$; and
\item for every $k\in [1,m]$, $\mathfrak{R}^{(2)}_{k}$ is a head-constant echelonless second-order path  that is not coherent,
\end{enumerate}
then
$\bb{\mathfrak{R}^{(2)}_{0}}<\bb{\mathfrak{R}^{(2)}_{m}}$.

Now, let $\mathbf{w}$ be a word in $S^{\star}$ of length $\bb{\mathbf{w}}=m+2$ and let $(\mathfrak{R}^{(2)}_{k})_{k\in\bb{\mathbf{w}}}$ is a family of second-order paths in $\mathrm{Pth}_{\boldsymbol{\mathcal{A}}^{(2)},\mathbf{w}}$ satisfying that 
\begin{enumerate}
\item for every $k\in m$, $(\mathfrak{R}^{(2)}_{k}, w_{k})\prec_{\mathbf{Pth}_{\boldsymbol{\mathcal{A}}^{(2)}}} (\mathfrak{R}^{(2)}_{k+1}, w_{k+1})$; and
\item for every $k\in [1,m]$, $\mathfrak{R}^{(2)}_{k}$ is a head-constant echelonless second-order path that is not coherent.
\end{enumerate}

Consider the final subsequence $(\mathfrak{R}^{(2)}_{k})_{k\in [1,m+1]}$. This is a sequence of length $m$ satisfying the conditions of the inductive step. Hence, $\bb{\mathfrak{R}^{(2)}_{1}}<\bb{\mathfrak{R}^{(2)}_{m+1}}$.

Now, consider the initial subsequence $(\mathfrak{R}^{(2)}_{k})_{k\in 2}$. This is a sequence of length $1$ instantiating that $(\mathfrak{R}^{(2)}_{0},w_{0})\prec_{\mathbf{Pth}_{\boldsymbol{\mathcal{A}}^{(2)}}}(\mathfrak{R}^{(2)}_{1},w_{1})$ and, by assumption, $\mathfrak{R}^{(2)}_{1}$ is a head-constant echelonless second-order path that is not coherent. By the base step, we have that $\bb{\mathfrak{R}^{(2)}_{0}}<\bb{\mathfrak{R}^{(2)}_{1}}$.

All in all, we can affirm that $\bb{\mathfrak{R}^{(2)}_{0}}<\bb{\mathfrak{R}^{(2)}_{m+1}}$.

This finishes the proof.
\end{proof}

In the following lemma we show that, for every pair of sorts $s$ and $t$ in $S$ and every pair of second-order paths $\mathfrak{P}^{(2)}\in\mathrm{Pth}_{\boldsymbol{\mathcal{A}}^{(2)},s}$ and $\mathfrak{Q}^{(2)}\in\mathrm{Pth}_{\boldsymbol{\mathcal{A}}^{(2)},t}$, if $(\mathfrak{Q}^{(2)},t)<_{\mathbf{Pth}_{\boldsymbol{\mathcal{A}}^{(2)}}}(\mathfrak{P}^{(2)},s)$ and there exists an strictly decreasing sequence instantiating this inequality composed of coherent, head-constant and echelonless second-order paths, then $\bb{\mathfrak{Q}^{(2)}}\leq \bb{\mathfrak{P}^{(2)}}$ and if $\mathfrak{Q}^{(2)}$ is a second-order path of length at least one, then the maximum of the heights of the first-order  translations occurring in $\mathfrak{Q}^{(2)}$ is strictly smaller than the maximum of the heights of the first-order  translations occurring in $\mathfrak{P}^{(2)}$.

\begin{lemma}\label{LDOrdV} 
Let $m$ be a natural number in $\mathbb{N}-\{0\}$, $\mathbf{w}$ a word in $S^{\star}$ such that $\bb{\mathbf{w}}=m+1$, and $(\mathfrak{R}^{(2)}_{k})_{k\in \bb{\mathbf{w}}}$ a family of second-order paths in $\mathrm{Pth}_{\boldsymbol{\mathcal{A}}^{(2)},\mathbf{w}}$ such that
\begin{enumerate}
\item for every $k\in m$, $(\mathfrak{R}^{(2)}_{k}, w_{k})\prec_{\mathbf{Pth}_{\boldsymbol{\mathcal{A}}^{(2)}}} (\mathfrak{R}^{(2)}_{k+1}, w_{k+1})$ and, 
\item for every $k\in [1,m]$, $\mathfrak{R}^{(2)}_{k}$ is a coherent head-constant echelonless second-order path.
\end{enumerate}
Then $\bb{\mathfrak{R}^{(2)}_{0}}\leq \bb{\mathfrak{R}^{(2)}_{m}}$. Moreover, if $\bb{\mathfrak{R}^{(2)}_{0}}\geq 1$, and the second-order paths $\mathfrak{R}^{(2)}_{0}$ and $\mathfrak{R}^{(2)}_{m}$ have the form
\begin{align*}
\mathfrak{R}^{(2)}_{0}&=\left(
([R_{0,j}]_{w_{0}})_{j\in\bb{\mathbf{c}}+1},
(\mathfrak{r}^{(2)}_{0,j})_{j\in\bb{\mathbf{c}}},
(T^{(1)}_{0,j})_{j\in\bb{\mathbf{c}}}
\right)
\\
\mathfrak{R}^{(2)}_{m}&=\left(
([R_{m,j}]_{w_{m}})_{j\in\bb{\mathbf{d}}+1},
(\mathfrak{r}^{(2)}_{m,j})_{j\in\bb{\mathbf{d}}},
(T^{(1)}_{m,j})_{j\in\bb{\mathbf{d}}}
\right),
\end{align*}
for some non-empty words $\mathbf{c},\mathbf{d}\in S^{\star}-\{\lambda\}$, then 
\[
\max\{
\bb{T^{(1)}_{0,j}}\mid j\in \bb{\mathbf{c}}
\}
<
\max\{
\bb{T^{(1)}_{m,j}}\mid j\in \bb{\mathbf{d}}
\}.
\]
\end{lemma}
\begin{proof}
We prove it by induction on $m\in\mathbb{N}-\{0\}$.

\textsf{Base step of the induction.}

For $m=1$, we have that $(\mathfrak{R}^{(2)}_{0},w_{0})\prec_{\mathbf{Pth}_{\boldsymbol{\mathcal{A}}^{(2)}}} (\mathfrak{R}^{(2)}_{1},w_{1})$, where $\mathfrak{R}^{(2)}_{1}$ is a coherent head-constant echelonless second-order $\mathbf{d}$-path in $\mathrm{Pth}_{\boldsymbol{\mathcal{A}}^{(2)},w_{1}}$, for some non-empty word $\mathbf{d}\in S^{\star}-\{\lambda\}$, of the form
$$
\mathfrak{R}^{(2)}_{1}=
\left(
([R_{1,i}]_{w_{1}})_{i\in\bb{\mathbf{d}}+1},
(\mathfrak{p}^{(2)}_{1,i})_{i\in\bb{\mathbf{d}}},
(T^{(1)}_{1,i})_{i\in\bb{\mathbf{d}}}
\right).
$$
such that, for a unique word $\mathbf{w_{1}}\in S^{\star}-\{\lambda\}$, a unique operation symbol $\tau\in\Sigma^{\boldsymbol{\mathcal{A}}}_{\mathbf{w_{1}},w_{1}}$, 
$(T^{(1)}_{1,i})_{i\in\bb{\mathbf{d}}}$ is a family of first-order  translations of type $\tau$.

Let $((\mathbf{d}_{j})_{j\in\bb{\mathbf{w_{1}}}}, (\mathfrak{R}^{(2)}_{1,j})_{j\in\bb{\mathbf{w_{1}}}})$ be the pair in $(S^{\star})^{\bb{\mathbf{w_{1}}}}\times \mathrm{Pth}_{\boldsymbol{\mathcal{A}}^{(2)},\mathbf{w_{1}}}$ we can extract from $\mathfrak{R}^{(2)}_{1}$ after applying the extraction procedure from Lemma~\ref{LDPthExtract}. Since $(\mathfrak{R}^{(2)}_{0},w_{0})\prec_{\mathbf{Pth}_{\boldsymbol{\mathcal{A}}^{(2)}}}(\mathfrak{R}^{(2)}_{1},w_{1})$, according to Definition~\ref{DDOrd}, there exists some $j\in\bb{\mathbf{w_{1}}}$ for which $w_{0}=\mathbf{w}_{\mathbf{1},j}$ and $\mathfrak{R}^{(2)}_{0}=\mathfrak{R}^{(2)}_{1,j}$. Let us recall from the proof of Lemma~\ref{LDPthExtract}, that $\sum_{j\in\bb{\mathbf{w_{1}}}}\bb{\mathbf{d}_{j}}=\bb{\mathbf{d}}$. Consequently, 
$$\bb{\mathfrak{R}^{(2)}_{0}}=\bb{\mathbf{d}_{j}}\leq\bb{\mathbf{d}}=\bb{\mathfrak{R}^{(2)}_{1}}.$$ 

Moreover, if $\bb{\mathfrak{R}^{(2)}_{0}}\geq 1$ and $\mathfrak{R}^{(2)}_{0}$ has the form
$$\mathfrak{R}^{(2)}_{0}=\left(
([R_{0,j}]_{w_{0}})_{j\in\bb{\mathbf{d}_{j}}+1},
(\mathfrak{r}^{(2)}_{0,j})_{j\in\bb{\mathbf{d}_{j}}},
(T^{(1)}_{0,j})_{j\in\bb{\mathbf{d}_{j}}}
\right),$$
we have, according to the proof of Lemma~\ref{LDPthExtract}, that $(T^{(1)}_{0,j})_{j\in\bb{\mathbf{d}_{j}}}$ is a family of derived first-order  translations of a subsequence of $(T^{(1)}_{1,i})_{i\in\bb{\mathbf{d}}}$, hence all the first-order  translations occurring in $\mathfrak{R}^{(2)}_{0}$ have their heights strictly bounded by the height of a first-order  translation occurring in $\mathfrak{R}^{(2)}_{1}$, hence
\[
\max\{
\bb{T^{(1)}_{0,j}}\mid j\in\bb{\mathbf{d}_{j}}
\}
<
\max\{
\bb{T^{(1)}_{1,i}}\mid i\in \bb{\mathbf{d}}
\}.
\]

\textsf{Inductive step of the induction.}

Let us suppose that the statement holds for sequences of length up to $m\in\mathbb{N}-\{0\}$. That is, if $\mathbf{w}$ is a word in $S^{\star}$ of length $\bb{\mathbf{w}}=m+1$ and  $(\mathfrak{R}^{(2)}_{k})_{k\in \bb{\mathbf{w}}}$ a family of second-order paths in $\mathrm{Pth}_{\boldsymbol{\mathcal{A}}^{(2)},\mathbf{w}}$ satisfying
\begin{enumerate}
\item for every $k\in m$, $(\mathfrak{R}^{(2)}_{k}, w_{k})\prec_{\mathbf{Pth}_{\boldsymbol{\mathcal{A}}^{(2)}}} (\mathfrak{R}^{(2)}_{k+1}, w_{k+1})$; and
\item for every $k\in [1,m]$, $\mathfrak{R}^{(2)}_{k}$ is a coherent head-constant echelonless second-order path,
\end{enumerate}
then $\bb{\mathfrak{R}^{(2)}_{0}}\leq \bb{\mathfrak{R}^{(2)}_{m}}$.

Moreover, if $\bb{\mathfrak{R}^{(2)}_{0}}\geq 1$, and the second-order paths $\mathfrak{R}^{(2)}_{0}$ and $\mathfrak{R}^{(2)}_{m}$ have the form
\begin{align*}
\mathfrak{R}^{(2)}_{0}&=\left(
([R_{0,j}]_{w_{0}})_{j\in\bb{\mathbf{c}}+1},
(\mathfrak{r}^{(2)}_{0,j})_{j\in\bb{\mathbf{c}}},
(T^{(1)}_{0,j})_{j\in\bb{\mathbf{c}}}
\right)
\\
\mathfrak{R}^{(2)}_{m}&=\left(
([R_{m,j}]_{w_{m}})_{j\in\bb{\mathbf{d}}+1},
(\mathfrak{r}^{(2)}_{m,j})_{j\in\bb{\mathbf{d}}},
(T^{(1)}_{m,j})_{j\in\bb{\mathbf{d}}}
\right),
\end{align*}
for some non-empty words $\mathbf{c},\mathbf{d}\in S^{\star}-\{\lambda\}$, then 
\[
\max\{
\bb{T^{(1)}_{0,j}}\mid j\in\bb{\mathbf{c}}
\}
<
\max\{
\bb{T^{(1)}_{m,j}}\mid j\in \bb{\mathbf{d}}
\}.
\]

Now, let $\mathbf{w}$ be a word in $S^{\star}$ of length $\bb{\mathbf{w}}=m+2$ and $(\mathfrak{R}^{(2)}_{k})_{k\in \bb{\mathbf{w}}}$ a family of second-order paths in $\mathrm{Pth}_{\boldsymbol{\mathcal{A}}^{(2)},\mathbf{w}}$ satisfying
\begin{enumerate}
\item for every $k\in m+1$, $(\mathfrak{R}^{(2)}_{k}, w_{k})\prec_{\mathbf{Pth}_{\boldsymbol{\mathcal{A}}^{(2)}}} (\mathfrak{R}^{(2)}_{k+1}, w_{k+1})$; and
\item for every $k\in [1,m+1]$, $\mathfrak{R}^{(2)}_{k}$ is a coherent head-constant echelonless second-order path.
\end{enumerate}

Consider the final subsequence $(\mathfrak{R}^{(2)}_{k})_{k\in[1,m+1]}$. This is a sequence of length $m$ satisfying the conditions of the inductive step. Hence, 
$\bb{\mathfrak{R}^{(2)}_{1}}\leq\bb{\mathfrak{R}^{(2)}_{m+1}}$. 

Moreover, by assumption, $\mathfrak{R}^{(2)}_{1}$ is a coherent head-constant echelonless second-order path. Thus, if $\mathfrak{R}^{(2)}_{1}$ and $\mathfrak{R}^{(2)}_{m+1}$ have the form
\begin{align*}
\mathfrak{R}^{(2)}_{1}&=\left(
([R_{1,j}]_{w_{1}})_{j\in\bb{\mathbf{d}}+1},
(\mathfrak{r}^{(2)}_{1,j})_{j\in\bb{\mathbf{d}}},
(T^{(1)}_{1,j})_{j\in\bb{\mathbf{d}}}
\right)
\\
\mathfrak{R}^{(2)}_{m+1}&=\left(
([R_{m+1,j}]_{w_{m+1}})_{j\in\bb{\mathbf{e}}+1},
(\mathfrak{r}^{(2)}_{m+1,j})_{j\in\bb{\mathbf{e}}},
(T^{(1)}_{m+1,j})_{j\in\bb{\mathbf{e}}}
\right),
\end{align*}
for some non-empty words $\mathbf{d},\mathbf{e}\in S^{\star}-\{\lambda\}$, then 
\[
\max\{
\bb{T^{(1)}_{1,j}}\mid j\in\bb{\mathbf{d}}
\}
<
\max\{
\bb{T^{(1)}_{m+1,j}}\mid j\in \bb{\mathbf{e}}
\}.
\]

Now, consider the initial subsequence $(\mathfrak{R}^{(2)}_{k})_{k\in 2}$. This is a sequence of length $1$ instantiating that $(\mathfrak{R}^{(2)}_{0},w_{0})\prec_{\mathbf{Pth}_{\boldsymbol{\mathcal{A}}^{(2)}}}(\mathfrak{R}^{(2)}_{1},w_{1})$ and, by assumption, $\mathfrak{R}^{(2)}_{1}$ is a coherent head-constant echelonless second-order path. By the base step, we have that $\bb{\mathfrak{R}^{(2)}_{0}}\leq \bb{\mathfrak{R}^{(2)}_{1}}$.

All in all, we can affirm that $\bb{\mathfrak{R}^{(2)}_{0}}\leq\bb{\mathfrak{R}^{(2)}_{m+1}}$.

Moreover, if $\bb{\mathfrak{R}^{(2)}_{0}}\geq 1$ and the second-order path $\mathfrak{R}^{(2)}_{0}$ has the form
\begin{align*}
\mathfrak{R}^{(2)}_{0}&=\left(
([R_{0,j}]_{w_{0}})_{j\in\bb{\mathbf{c}}+1},
(\mathfrak{r}^{(2)}_{0,j})_{j\in\bb{\mathbf{c}}},
(T^{(1)}_{0,j})_{j\in\bb{\mathbf{c}}}
\right),
\end{align*}
for some non-empty word $\mathbf{c}\in S^{\star}-\{\lambda\}$, then 
\[
\max\{
\bb{T^{(1)}_{0,j}}\mid j\in\bb{\mathbf{c}}
\}
<
\max\{
\bb{T^{(1)}_{1,j}}\mid j\in\bb{\mathbf{d}}
\}.
\]

All in all, we can affirm in this case that 
\[
\max\{
\bb{T^{(1)}_{0,j}}\mid j\in\bb{\mathbf{c}}
\}
<
\max\{
\bb{T^{(1)}_{m+1,j}}\mid j\in \bb{\mathbf{e}}
\}.
\]

This finishes the proof.
\end{proof}

In the following corollary we show that, for every pair of sorts $s$ and $t$ in $S$ and every pair of second-order paths $\mathfrak{P}^{(2)}\in\mathrm{Pth}_{\boldsymbol{\mathcal{A}}^{(2)},s}$ and $\mathfrak{Q}^{(2)}\in\mathrm{Pth}_{\boldsymbol{\mathcal{A}}^{(2)},t}$, if $(\mathfrak{Q}^{(2)},t)<_{\mathbf{Pth}_{\boldsymbol{\mathcal{A}}^{(2)}}}(\mathfrak{P}^{(2)},s)$ and $\mathfrak{P}^{(2)}$ is either (1) a second-order path of length strictly greater than one containing at least one second-order echelon; or (2) an echelonless second-order path that is not head-constant; or (3) a head-constant echelonless second-order path that is not coherent, then $\bb{\mathfrak{Q}^{(2)}}<\bb{\mathfrak{P}^{(2)}}$.

\begin{corollary}\label{CDOrdI}Let $m$ be a natural number in $\mathbb{N}-\{0\}$, $\mathbf{w}$ a word in $S^{\star}$ such that $\bb{\mathbf{w}}=m+1$, and $(\mathfrak{R}^{(2)}_{k})_{k\in \bb{\mathbf{w}}}$ a family of second-order paths in $\mathrm{Pth}_{\boldsymbol{\mathcal{A}}^{(2)},\mathbf{w}}$ such that
\begin{enumerate}
\item  for every $k\in m$, $(\mathfrak{R}^{(2)}_{k}, w_{k})\prec_{\mathbf{Pth}_{\boldsymbol{\mathcal{A}}^{(2)}}} (\mathfrak{R}^{(2)}_{k+1}, w_{k+1})$;
\end{enumerate}
and one of the following conditions is met
\begin{enumerate}
\item[(2.1)]$\mathfrak{R}^{(2)}_{m}$ is a second-order path of length strictly greater than one containing at least one  second-order echelon; or 
\item[(2.2)]$\mathfrak{R}^{(2)}_{m}$ an echelonless second-order path  that is not head-constant; or 
\item[(2.3)] $\mathfrak{R}^{(2)}_{m}$ a head-constant echelonless second-order path  that is not coherent.
\end{enumerate}
Then $\bb{\mathfrak{R}^{(2)}_{0}}<\bb{\mathfrak{R}^{(2)}_{m}}$.
\end{corollary}

\begin{proof}
We prove it by induction on $m\in\mathbb{N}-\{0\}$.

\textsf{Base step of the induction.}

For $m=1$, we have that $(\mathfrak{R}^{(2)}_{0},w_{0})\prec_{\mathbf{Pth}_{\boldsymbol{\mathcal{A}}^{(2)}}}(\mathfrak{R}^{(2)}_{1},w_{1})$, where either (2.1) $\mathfrak{R}^{(2)}_{1}$ is a second-order path of length strictly greater than one containing at least one   second-order echelon; or 
(2.2) $\mathfrak{R}^{(2)}_{1}$ an echelonless second-order path  that is not head-constant; or 
(2.3) $\mathfrak{R}^{(2)}_{1}$ a head-constant echelonless second-order path that is not coherent. Then by either Lemma~\ref{LDOrdII},~\ref{LDOrdIII}, or~\ref{LDOrdIV}, respectively, we have that 
$\bb{\mathfrak{R}^{(2)}_{0}}<\bb{\mathfrak{R}^{(2)}_{1}}$.

\textsf{Inductive step of the induction.}

Let us suppose that the statement holds for sequences of length up to $m\in\mathbb{N}-\{0\}$. That is, if $\mathbf{w}$ is a word in $S^{\star}$ of length $\bb{\mathbf{w}}=m+1$ and $(\mathfrak{R}^{(2)}_{k})_{k\in \bb{\mathbf{w}}}$ a family of second-order paths in $\mathrm{Pth}_{\boldsymbol{\mathcal{A}}^{(2)},\mathbf{w}}$ such that
\begin{enumerate}
\item  for every $k\in m$, $(\mathfrak{R}^{(2)}_{k}, w_{k})\prec_{\mathbf{Pth}_{\boldsymbol{\mathcal{A}}^{(2)}}} (\mathfrak{R}^{(2)}_{k+1}, w_{k+1})$;
\end{enumerate}
and either one of the following conditions holds
\begin{enumerate}
\item[(2.1)]$\mathfrak{R}^{(2)}_{m}$ is a second-order path of length strictly greater than one containing at least one   second-order echelon; or 
\item[(2.2)]$\mathfrak{R}^{(2)}_{m}$ an echelonless second-order path  that is not head-constant; or 
\item[(2.3)] $\mathfrak{R}^{(2)}_{m}$ a head-constant echelonless second-order path that is not coherent;
\end{enumerate}
then $\bb{\mathfrak{R}^{(2)}_{0}}<\bb{\mathfrak{R}^{(2)}_{m}}$.

Now, let $\mathbf{w}$ be a word in $S^{\star}$ of length $\bb{\mathbf{w}}=m+2$ and let $(\mathfrak{R}^{(2)}_{k})_{k\in \bb{\mathbf{w}}}$ be a family of second-order paths in $\mathrm{Pth}_{\boldsymbol{\mathcal{A}}^{(2)},\mathbf{w}}$ such that
\begin{enumerate}
\item  for every $k\in m$, $(\mathfrak{R}^{(2)}_{k}, w_{k})\prec_{\mathbf{Pth}_{\boldsymbol{\mathcal{A}}^{(2)}}} (\mathfrak{R}^{(2)}_{k+1}, w_{k+1})$;
\end{enumerate}
and either one of the following conditions holds
\begin{enumerate}
\item[(2.1)]$\mathfrak{R}^{(2)}_{m}$ is a second-order path of length strictly greater than one containing at least one   second-order echelon; or 
\item[(2.2)]$\mathfrak{R}^{(2)}_{m}$ an echelonless second-order path  that is not head-constant; or 
\item[(2.3)] $\mathfrak{R}^{(2)}_{m}$ a head-constant echelonless second-order path  that is not coherent.
\end{enumerate}

Consider the final subsequence $(\mathfrak{R}^{(2)}_{k})_{k\in[1,m+1]}$. This is a sequence of length $m$ satisfying the conditions of the inductive step. Hence $\bb{\mathfrak{R}^{(2)}_{1}}<\bb{\mathfrak{R}^{(2)}_{m+1}}$.

Now, consider the initial subsequence $(\mathfrak{R}^{(2)}_{k})_{k\in 2}$. This is a sequence of length $1$ instantiating that $(\mathfrak{R}^{(2)}_{0},w_{0})\prec_{\mathbf{Pth}_{\boldsymbol{\mathcal{A}}^{(2)}}} (\mathfrak{R}^{(2)}_{1},w_{1})$. Since $\prec_{\mathbf{Pth}_{\boldsymbol{\mathcal{A}}^{(2)}}}$ is an irreflexive relation by Remark~\ref{RDIrrefl}, we conclude that $(\mathfrak{R}^{(2)}_{1},w_{1})$ cannot be 
a minimal element of $\leq_{\mathbf{Pth}_{\boldsymbol{\mathcal{A}}^{(2)}}}$. Thus, taking into account Lemma~\ref{LDOrdI}, the different possibilities for $\mathfrak{R}^{(2)}_{1}$ are
\begin{enumerate}
\item[(1)] $\mathfrak{R}^{(2)}_{1}$ is a $(2,[1])$-identity second-order path on a non-minimal path term class. Following Proposition~\ref{PDULower}, we conclude that $\mathfrak{R}^{(2)}_{0}$ is also a $(2,[1])$-identity second-order path. Hence, $\bb{\mathfrak{R}^{(2)}_{0}}=\bb{\mathfrak{R}^{(2)}_{1}}$; or
\item[(2)] $\mathfrak{R}^{(2)}_{1}$ is a second-order path of length striclty greater than one containing at least one   second-order echelon. By Lemma~\ref{LDOrdII} we conclude that $\bb{\mathfrak{R}^{(2)}_{0}}<\bb{\mathfrak{R}^{(2)}_{1}}$; or
\item[(3)] $\mathfrak{R}^{(2)}_{1}$ is an echelonless second-order path that is not head-constant. By Lemma~\ref{LDOrdIII} we conclude that $\bb{\mathfrak{R}^{(2)}_{0}}<\bb{\mathfrak{R}^{(2)}_{1}}$; or
\item[(4)] $\mathfrak{R}^{(2)}_{1}$ is a head-constant echelonless second-order path that is not coherent. By Lemma~\ref{LDOrdIV} we conclude that $\bb{\mathfrak{R}^{(2)}_{0}}<\bb{\mathfrak{R}^{(2)}_{1}}$; or
\item[(5)] $\mathfrak{R}^{(2)}_{1}$ is a coherent head-constant echelonless second-order path. By Lemma~\ref{LDOrdV} we conclude that $\bb{\mathfrak{R}^{(2)}_{0}}\leq\bb{\mathfrak{R}^{(2)}_{1}}$.
\end{enumerate}
All in all, we conclude that 
$\bb{\mathfrak{R}^{(2)}_{0}}<\bb{\mathfrak{R}^{(2)}_{m+1}}$.

This finishes the proof.
\end{proof}

In the following corollary we show that, for every pair of sorts $s$ and $t$ in $S$ and every pair of second-order paths $\mathfrak{P}^{(2)}\in\mathrm{Pth}_{\boldsymbol{\mathcal{A}}^{(2)},s}$ and $\mathfrak{Q}^{(2)}\in\mathrm{Pth}_{\boldsymbol{\mathcal{A}}^{(2)},t}$, if $(\mathfrak{Q}^{(2)},t)<_{\mathbf{Pth}_{\boldsymbol{\mathcal{A}}^{(2)}}}(\mathfrak{P}^{(2)},s)$ and $\mathfrak{P}^{(2)}$ is either (1) a second-order path of length strictly greater than one containing at least one   second-order echelon; or (2) an echelonless second-order path, then $\bb{\mathfrak{Q}^{(2)}}\leq \bb{\mathfrak{P}^{(2)}}$.

\begin{corollary}\label{CDOrdII}Let $m$ be a natural number in $\mathbb{N}-\{0\}$, $\mathbf{w}$ a word in $S^{\star}$ such that $\bb{\mathbf{w}}=m+1$, and  $(\mathfrak{R}^{(2)}_{k})_{k\in \bb{\mathbf{w}}}$ a family of second-order paths in $\mathrm{Pth}_{\boldsymbol{\mathcal{A}}^{(2)},\mathbf{w}}$ such that
\begin{enumerate}
\item for every $k\in m$, $(\mathfrak{R}^{(2)}_{k}, w_{k})\prec_{\mathbf{Pth}_{\boldsymbol{\mathcal{A}}^{(2)}}} (\mathfrak{R}^{(2)}_{k+1}, w_{k+1})$;
\end{enumerate}
and one of the following conditions is met
\begin{enumerate}
\item[(2.1)] $\mathfrak{R}^{(2)}_{m}$ is a second-order path of length strictly greater than one containing at least one   second-order echelon; or 
\item[(2.2)] $\mathfrak{R}^{(2)}_{m}$ an echelonless second-order path.
\end{enumerate}
Then $\bb{\mathfrak{R}^{(2)}_{0}}\leq \bb{\mathfrak{R}^{(2)}_{m}}$.
\end{corollary}
\begin{proof}
We prove it by induction on $m\in\mathbb{N}-\{0\}$.

\textsf{Base step of the induction.}

For $m=1$, we have that  $(\mathfrak{R}^{(2)}_{0},w_{0})\prec_{\mathbf{Pth}_{\boldsymbol{\mathcal{A}}^{(2)}}}(\mathfrak{R}^{(2)}_{1},w_{1})$, where either (2.1) $\mathfrak{R}^{(2)}_{1}$ is a second-order path of length strictly greater than one containing at least one   second-order echelon; or  $\mathfrak{R}^{(2)}_{1}$ an echelonless second-order path. Then by either Lemma~\ref{LDOrdII} or Lemmas~\ref{LDOrdIII},~\ref{LDOrdIV} or~\ref{LDOrdV}, respectively, we have that $\bb{\mathfrak{R}^{(2)}_{0}}\leq\bb{\mathfrak{R}^{(2)}_{1}}$.

\textsf{Inductive step of the induction.}

Let us suppose that the statement holds for sequences of length up to $m\in\mathbb{N}-\{0\}$. That is, if $\mathbf{w}$ is a word in $S^{\star}$ of length $\bb{\mathbf{w}}=m+1$ and $(\mathfrak{R}^{(2)}_{k})_{k\in\bb{\mathbf{w}}}$ a family of second-order paths in $\mathrm{Pth}_{\boldsymbol{\mathcal{A}}^{(2)},\mathbf{w}}$ such that
\begin{enumerate}
\item for every $k\in m$, $(\mathfrak{R}^{(2)}_{k}, w_{k})\prec_{\mathbf{Pth}_{\boldsymbol{\mathcal{A}}^{(2)}}} (\mathfrak{R}^{(2)}_{k+1}, w_{k+1})$;
\end{enumerate}
and  either one of the following conditions holds
\begin{enumerate}
\item[(2.1)] $\mathfrak{R}^{(2)}_{m}$ is a second-order path of length strictly greater than one containing at least one   second-order echelon; or 
\item[(2.2)] $\mathfrak{R}^{(2)}_{m}$ an echelonless second-order path;
\end{enumerate}
then $\bb{\mathfrak{R}^{(2)}_{0}}\leq \bb{\mathfrak{R}^{(2)}_{m}}$.

Now let $\mathbf{w}$ be a word in $S^{\star}$ of length $\bb{\mathbf{w}}=m+2$ and let $(\mathfrak{R}^{(2)}_{k})_{k\in\bb{\mathbf{w}}}$ be a family of second-order paths in $\mathrm{Pth}_{\boldsymbol{\mathcal{A}}^{(2)},\mathbf{w}}$ such that
\begin{enumerate}
\item for every $k\in m+1$, $(\mathfrak{R}^{(2)}_{k}, w_{k})\prec_{\mathbf{Pth}_{\boldsymbol{\mathcal{A}}^{(2)}}} (\mathfrak{R}^{(2)}_{k+1}, w_{k+1})$;
\end{enumerate}
and  either one of the following conditions holds
\begin{enumerate}
\item[(2.1)] $\mathfrak{R}^{(2)}_{m+1}$ is a second-order path of length strictly greater than one containing at least one   second-order echelon; or 
\item[(2.2)] $\mathfrak{R}^{(2)}_{m+1}$ an echelonless second-order path.
\end{enumerate}

Consider the final subsequence $(\mathfrak{R}^{(2)}_{k})_{k\in[1,m+1]}$. This is a sequence of length $m$ satisfying the conditions of the inductive step. Hence $\bb{\mathfrak{R}^{(2)}_{1}}\leq\bb{\mathfrak{R}^{(2)}_{m+1}}$.

Now, consider the initial subsequence $(\mathfrak{R}^{(2)}_{k})_{k\in 2}$. This is a sequence of length $1$ instantiating that $(\mathfrak{R}^{(2)}_{0},w_{0})\prec_{\mathbf{Pth}_{\boldsymbol{\mathcal{A}}^{(2)}}} (\mathfrak{R}^{(2)}_{1},w_{1})$. Since $\prec_{\mathbf{Pth}_{\boldsymbol{\mathcal{A}}^{(2)}}}$ is an irreflexive relation by Remark~\ref{RDIrrefl}, we conclude that $(\mathfrak{R}^{(2)}_{1},w_{1})$ cannot be a minimal element of $\leq_{\mathbf{Pth}_{\boldsymbol{\mathcal{A}}^{(2)}}}$. Thus, taking into account Lemma~\ref{LDOrdI}, the different possibilities for $\mathfrak{R}^{(2)}_{1}$ are
\begin{enumerate}
\item[(1)] $\mathfrak{R}^{(2)}_{1}$ is a $(2,[1])$-identity second-order path on a non-minimal path term class. Following Proposition~\ref{PDULower}, we conclude that $\mathfrak{R}^{(2)}_{0}$ is also a $(2,[1])$-identity second-order path. Hence, $\bb{\mathfrak{R}^{(2)}_{0}}=\bb{\mathfrak{R}^{(2)}_{1}}$; or
\item[(2)] $\mathfrak{R}^{(2)}_{1}$ is a second-order path of length strictly greater than one containing at least one   second-order echelon. By Lemma~\ref{LDOrdII} we conclude that 
$\bb{\mathfrak{R}^{(2)}_{0}}<\bb{\mathfrak{R}^{(2)}_{1}}$; or
\item[(3)] $\mathfrak{R}^{(2)}_{1}$ is an echelonless second-order path. Then by Lemmas~\ref{LDOrdIII},~\ref{LDOrdIV} or~\ref{LDOrdV}, we have that $\bb{\mathfrak{R}^{(2)}_{0}}\leq\bb{\mathfrak{R}^{(2)}_{1}}$.
\end{enumerate}
All in all, we conclude that $\bb{\mathfrak{R}^{(2)}_{0}}\leq\bb{\mathfrak{R}^{(2)}_{m+1}}$.

This finishes the proof.
\end{proof}

The above technical lemmas will be used in the proof of the following proposition.

\begin{restatable}{proposition}{PDOrdArt}
\label{PDOrdArt}
The preorder $\leq_{\mathbf{Pth}_{\boldsymbol{\mathcal{A}}^{(2)}}}$ on $\coprod\mathrm{Pth}_{\boldsymbol{\mathcal{A}}^{(2)}}$ is antisymmetric and, then, in the ordered set $(\coprod\mathrm{Pth}_{\boldsymbol{\mathcal{A}}^{(2)}}, \leq_{\mathbf{Pth}_{\boldsymbol{\mathcal{A}}^{(2)}}})$ there is not any strictly decreasing $\omega_{0}$-chain, i.e., $(\coprod\mathrm{Pth}_{\boldsymbol{\mathcal{A}}^{(2)}}, \leq_{\mathbf{Pth}_{\boldsymbol{\mathcal{A}}^{(2)}}})$ is an Artinian ordered set.
\end{restatable}
\begin{proof}
We first prove that $\leq_{\mathbf{Pth}_{\boldsymbol{\mathcal{A}}^{(2)}}}$ is antisymmetric. 

Let $s$ and $t$ be sorts in $S$, $\mathfrak{P}^{(2)}\in\mathrm{Pth}_{\boldsymbol{\mathcal{A}}^{(2)},s}$, and $\mathfrak{Q}^{(2)}\in\mathrm{Pth}_{\boldsymbol{\mathcal{A}}^{(2)},t}$. Let us suppose that $(\mathfrak{Q}^{(2)},t)\leq_{\mathbf{Pth}_{\boldsymbol{\mathcal{A}}^{(2)}}} (\mathfrak{P}^{(2)},s)$ and $(\mathfrak{P}^{(2)},s)\leq_{\mathbf{Pth}_{\boldsymbol{\mathcal{A}}^{(2)}}} (\mathfrak{Q}^{(2)},t)$. We want to show that $s=t$ and $\mathfrak{Q}^{(2)}=\mathfrak{P}^{(2)}$.

This statement is trivial if $\mathfrak{P}^{(2)}$ or $\mathfrak{Q}^{(2)}$ is minimal in $(\coprod\mathrm{Pth}_{\boldsymbol{\mathcal{A}}^{(2)}},\leq_{\mathbf{Pth}_{\boldsymbol{\mathcal{A}}^{(2)}}})$. So let us assume that neither $\mathfrak{P}^{(2)}$ nor $\mathfrak{Q}^{(2)}$ are minimal in $(\coprod\mathrm{Pth}_{\boldsymbol{\mathcal{A}}^{(2)}},\leq_{\mathbf{Pth}_{\boldsymbol{\mathcal{A}}^{(2)}}})$. Following Lemma~\ref{LDOrdI}, for $\mathfrak{P}^{(2)}$ we have that either (A) it is a $(2,[1])$-identity second-order path on a non-minimal path term class, or (B) it is a second-order path of length strictly greater than one containing at least one   second-order echelon or (C) it is an echelonless second-order path. In the same way, for $\mathfrak{Q}^{(2)}$ we have that either (D) it is a $(2,[1])$-identity second-order path on a non-minimal path term class, or (E) it is a second-order path of length strictly greater than one containing at least one   second-order echelon or (F) it is an echelonless second-order path. However, as we will show immediately below only cases $(A,D)$, $(B,E)$ and $(C,F)$ are feasible.

We first consider the case where (A) $\mathfrak{P}^{(2)}$ is a $(2,[1])$-identity second-order path. Since we are assuming that $(\mathfrak{Q}^{(2)},t)\leq_{\mathbf{Pth}_{\boldsymbol{\mathcal{A}}^{(2)}}} (\mathfrak{P}^{(2)},s)$ we conclude, in virtue of Proposition~\ref{PDULower} that (D) $\mathfrak{Q}^{(2)}$ is a $(2,[1])$-identity second-order path. Thus, if the $(2,[1])$-identity second-order paths $\mathfrak{P}^{(2)}$ and $\mathfrak{Q}^{(2)}$ have the form
\begin{align*}
\mathfrak{Q}^{(2)}&=\mathrm{ip}^{(2,[1])\sharp}_{t}\left([Q]_{t}
\right),
&
\mathfrak{P}^{(2)}&=\mathrm{ip}^{(2,[1])\sharp}_{s}\left(
[P]_{s}
\right),
\end{align*}
for a suitable pair of path term classes $[Q]_{t}$ in $[\mathrm{PT}_{\boldsymbol{\mathcal{A}}}]_{t}$ and $[P]_{s}$ in $[\mathrm{PT}_{\boldsymbol{\mathcal{A}}}]_{s}$. Now, since we are assuming that $(\mathfrak{Q}^{(2)},t)\leq_{\mathbf{Pth}_{\boldsymbol{\mathcal{A}}^{(2)}}} (\mathfrak{P}^{(2)},s)$ and $(\mathfrak{P}^{(2)},s)\leq_{\mathbf{Pth}_{\boldsymbol{\mathcal{A}}^{(2)}}} (\mathfrak{Q}^{(2)},t)$, we conclude that 
$([Q]_{t},t)\leq_{[\mathbf{PT}_{\boldsymbol{\mathcal{A}}}]} ([P]_{s},s)$ and $([P]_{s},s)\leq_{[\mathbf{PT}_{\boldsymbol{\mathcal{A}}}]} ([Q]_{t},t)$. 

Since $(\coprod[\mathrm{PT}_{\boldsymbol{\mathcal{A}}}], \leq_{[\mathbf{PT}_{\boldsymbol{\mathcal{A}}}]})$ is an order according to Proposition~\ref{PPTQOrdArt}, we conclude that $([Q]_{t},t)=([P]_{s},s)$. In particular $s=t$ and 
$$\mathfrak{Q}^{(2)}=\mathrm{ip}^{(2,[1])\sharp}_{t}
\left([Q]_{t}
\right)=\mathrm{ip}^{(2,[1])\sharp}_{s}\left(
[P]_{s}
\right)=\mathfrak{P}^{(2)}.$$

The same result can be obtained by assuming that (D) $\mathfrak{Q}^{(2)}$ is a $(2,[1])$-identity second-order path. Therefore cases (A) and (D) are totally determined.

We next show that cases (B,F) and (C,E) are infeasible.

If (B,F), then $\mathfrak{P}^{(2)}$ and $\mathfrak{Q}^{(2)}$ are necessarily different since one contains at least one   second-order echelon and the other is echelonless.  But $(\mathfrak{Q}^{(2)},t)<_{\mathbf{Pth}_{\boldsymbol{\mathcal{A}}^{(2)}}} (\mathfrak{P}^{(2)},s)$ and $\mathfrak{P}^{(2)}$ is a second-order path of length strictly greater than one containing at least one   second-order echelon, hence, by Corollary~\ref{CDOrdI}, we have that $\bb{\mathfrak{Q}^{(2)}}<\bb{\mathfrak{P}^{(2)}}$. Moreover, since $(\mathfrak{P}^{(2)},s)<_{\mathbf{Pth}_{\boldsymbol{\mathcal{A}}^{(2)}}} (\mathfrak{Q}^{(2)},t)$ and $\mathfrak{Q}^{(2)}$ is an echelonless second-order path, we have, by Corollary~\ref{CDOrdII}, that $\bb{\mathfrak{P}^{(2)}}\leq \bb{\mathfrak{Q}^{(2)}}$. This leads to the following contradiction
$$
\bb{\mathfrak{P}^{(2)}}\leq\bb{\mathfrak{Q}^{(2)}}<\bb{\mathfrak{P}^{(2)}}.
$$

The same conclusion can be drawn for case (C,E), by applying to it the same reasoning as in the case (B,F).

Therefore we are only left with these options: (B,E) and (C,F), where the second-order paths $\mathfrak{P}^{(2)}$ and $\mathfrak{Q}^{(2)}$ are required to have the same nature.

For (B,E), let us suppose, towards a contradiction, that $\mathfrak{P}^{(2)}\neq\mathfrak{Q}^{(2)}$. Then, since $(\mathfrak{Q}^{(2)},t)<_{\mathbf{Pth}_{\boldsymbol{\mathcal{A}}^{(2)}}}(\mathfrak{P}^{(2)},s)$ and $\mathfrak{P}^{(2)}$ is a second-order path of length strictly greater than one containing at least one   second-order echelon, hence, by Corollary~\ref{CDOrdI}, we have that $\bb{\mathfrak{Q}^{(2)}}<\bb{\mathfrak{P}^{(2)}}$. Moreover, since $(\mathfrak{P}^{(2)},s)<_{\mathbf{Pth}_{\boldsymbol{\mathcal{A}}^{(2)}}}(\mathfrak{Q}^{(2)},t)$ and $\mathfrak{Q}^{(2)}$  is a second-order path of length strictly greater than one containing at least one   second-order echelon, hence, by Corollary~\ref{CDOrdI}, we have that $\bb{\mathfrak{P}^{(2)}}<\bb{\mathfrak{Q}^{(2)}}$. This leads to the following contradiction
$$
\bb{\mathfrak{P}^{(2)}}<\bb{\mathfrak{Q}^{(2)}}<\bb{\mathfrak{P}^{(2)}}.
$$
Therefore, we infer that $\mathfrak{P}^{(2)}=\mathfrak{Q}^{(2)}$.

For (C,F), let us suppose, towards a contradiction that $\mathfrak{P}^{(2)}\neq \mathfrak{Q}^{(2)}$. Then, since $(\mathfrak{Q}^{(2)},t)<_{\mathbf{Pth}_{\boldsymbol{\mathcal{A}}^{(2)}}} (\mathfrak{P}^{(2)},s)$, we have, by Remark~\ref{RDOrd}, that there exists a natural number $m\in\mathbb{N}-\{0\}$, a word $\mathbf{w}\in S^{\star}$ of length $\bb{\mathbf{w}}=m+1$, and a family of second-order paths $(\mathfrak{R}^{(2)}_{k})_{k\in\bb{\mathbf{w}}}$ in $\mathrm{Pth}_{\boldsymbol{\mathcal{A}}^{(2)},\mathbf{w}}$, such that $w_{0}=t$, $\mathfrak{R}^{(2)}_{0}=\mathfrak{Q}^{(2)}$, $w_{m}=s$, $\mathfrak{R}^{(2)}_{m}=\mathfrak{P}^{(2)}$ and, for every $k\in m$, $(\mathfrak{R}^{(2)}_{k},w_{k})\prec_{\mathbf{Pth}_{\boldsymbol{\mathcal{A}}^{(2)}}} (\mathfrak{R}^{(2)}_{k+1},w_{k+1})$. Then, by Lemma~\ref{LDOrdI}, for every $k\in [1,m]$, we have that 
\begin{enumerate}
\item $\mathfrak{R}^{(2)}_{k}$ is a $(2,[1])$-identity second-order path on a non-minimal path term class;
\item $\mathfrak{R}^{(2)}_{k}$ is a second-order path of length strictly greater than one containing at least one   second-order echelon; or 
\item $\mathfrak{R}^{(2)}_{k}$ is an echelonless second-order path.
\end{enumerate}

We claim that, among all the above possibilities, for every $k\in [1,m]$, the second-order path $\mathfrak{R}^{(2)}_{k}$ has to be a coherent head-constant echelonless second-order path. 

Let us suppose, towards a contradiction that there exists an index $k\in [1,m]$ for which $\mathfrak{R}^{(2)}_{k}$ is a $(2,[1])$-identity second-order path on a non-minimal path term class. Then, taking into account that $(\mathfrak{Q}^{(2)},t)\leq_{\mathbf{Pth}_{\boldsymbol{\mathcal{A}}^{(2)}}} (\mathfrak{R}^{(2)}_{k},w_{k})$, we obtain, 
according to Proposition~\ref{PDULower}, that $\mathfrak{Q}^{(2)}$ has to be a $(2,[1])$-identity second-order path, contradicting the nature of $\mathfrak{Q}^{(2)}$. 

Let us suppose, towards a contradiction that there exists an index $k\in [1,m]$ for which 
$\mathfrak{R}^{(2)}_{k}$ is a second-order path of length strictly greater than one containing at least one   second-order echelon. Then, taking into account that $(\mathfrak{R}^{(2)}_{k},w_{k})\leq_{\mathbf{Pth}_{\boldsymbol{\mathcal{A}}^{(2)}}} (\mathfrak{P}^{(2)},s)$ and $\mathfrak{P}^{(2)}$ is an echelonless second-order path, we obtain, according to Corollary~\ref{CDOrdII}, that 
$\bb{\mathfrak{R}^{(2)}_{k}}\leq\bb{\mathfrak{P}^{(2)}}$. Furthermore, since $(\mathfrak{Q}^{(2)},t)\leq_{\mathbf{Pth}_{\boldsymbol{\mathcal{A}}^{(2)}}} (\mathfrak{R}^{(2)}_{k},w_{k})$ we are assuming that $\mathfrak{R}^{(2)}_{k}$ is a second-order path of length strictly greater than one containing at least one   second-order echelon, we obtain, according to Corollary~\ref{CDOrdI}, that $\bb{\mathfrak{Q}^{(2)}}<\bb{\mathfrak{R}^{(2)}_{k}}$. Finally, since $(\mathfrak{P}^{(2)},s)\leq_{\mathbf{Pth}_{\boldsymbol{\mathcal{A}}^{(2)}}} (\mathfrak{Q}^{(2)},t)$ and we are working under the assumption that $\mathfrak{Q}^{(2)}$ is an echelonless second-order path, we obtain, according to Corollary~\ref{CDOrdII}, that $\bb{\mathfrak{P}^{(2)}}\leq \bb{\mathfrak{Q}^{(2)}}$. This leads to the following contradiction
$$
\bb{\mathfrak{P}^{(2)}}\leq\bb{\mathfrak{Q}^{(2)}}<\bb{\mathfrak{R}^{(2)}_{k}}\leq \bb{\mathfrak{P}^{(2)}}.
$$
Therefore, we infer that, for every $k\in [1,m]$, $\mathfrak{R}^{(2)}_{k}$ cannot be a second-order path of length strictly greater than one containing at least one   second-order echelon.

The same argument can be applied to obtain that, for every $k\in [1,m]$, $\mathfrak{R}^{(2)}_{k}$ cannot be an echelonless second-order path  that is not head-constant or a head-constant echelonless second-order path that is not coherent. 

All in all, we conclude that for every $k\in [1,m]$, $\mathfrak{R}^{(2)}_{k}$ is a coherent head-constant echelonless second-order path. In particular $\mathfrak{P}^{(2)}$ is a coherent head-constant echelonless second-order path.

Now, since $(\mathfrak{P}^{(2)},s)<_{\mathbf{Pth}_{\boldsymbol{\mathcal{A}}^{(2)}}} (\mathfrak{Q}^{(2)},t)$, we have, by Remark~\ref{RDOrd}, that there exists a natural number $n\in\mathbb{N}-\{0\}$, a word $\mathbf{v}\in S^{\star}$ of length $\bb{\mathbf{v}}=n+1$, and a family of second-order paths $(\mathfrak{S}^{(2)}_{l})_{l\in\bb{\mathbf{v}}}$ in $\mathrm{Pth}_{\boldsymbol{\mathcal{A}}^{(2)},\mathbf{v}}$, such that $v_{0}=s$, $\mathfrak{S}^{(2)}_{0}=\mathfrak{P}^{(2)}$, $v_{n}=t$, $\mathfrak{S}^{(2)}_{n}=\mathfrak{P}^{(2)}$ and, for every $l\in n$, $(\mathfrak{S}^{(2)}_{l},v_{l})\prec_{\mathbf{Pth}_{\boldsymbol{\mathcal{A}}^{(2)}}} (\mathfrak{S}^{(2)}_{l+1},v_{l+1})$.

By a proof similar to the one above, we can show that, for every $l\in [1,n]$, the second-order path $\mathfrak{S}^{(2)}_{l}$ has to be a coherent head-constant echelonless second-order path. In particular, $\mathfrak{Q}^{(2)}$ is also a coherent head-constant echelonless second-order path.

Note that this situation fulfils the premises of Lemma~\ref{LDOrdV}. 

Assume that, for a unique word $\mathbf{c}\in S^{\star}-\{\lambda\}$, $\mathfrak{P}^{(2)}$ is a second-order $\mathbf{c}$-path in $\boldsymbol{\mathcal{A}}^{(2)}$  of the form
$$
\mathfrak{P}^{(2)}
=
\left(
([P_{i}]_{s})_{i\in\bb{\mathbf{c}}+1},
(\mathfrak{p}^{(2)}_{i})_{i\in\bb{\mathbf{c}}},
(T^{(1)}_{i})_{i\in\bb{\mathbf{c}}}
\right),
$$
and, for a unique $\mathbf{d}\in S^{\star}$, $\mathfrak{Q}^{(2)}$ is a second-order $\mathbf{d}$-path in $\boldsymbol{\mathcal{A}}^{(2)}$  of the form
$$
\mathfrak{Q}^{(2)}
=
\left(
([Q_{j}]_{s})_{j\in\bb{\mathbf{d}}+1},
(\mathfrak{q}^{(2)}_{j})_{j\in\bb{\mathbf{d}}},
(U^{(1)}_{j})_{j\in\bb{\mathbf{d}}}
\right).
$$ 

Since $(\mathfrak{Q}^{(2)},t)<_{\mathbf{Pth}_{\boldsymbol{\mathcal{A}}^{(2)}}} (\mathfrak{P}^{(2)},s)$ and the sequence witnessing this fact is entirely composed of coherent head-constant echelonless second-order paths we have, by Lemma~\ref{LDOrdV}, that 
\[
\max\{
\bb{U^{(1)}_{j}}\mid j\in \bb{\mathbf{d}}
\}
<
\max\{
\bb{T^{(1)}_{i}}\mid i\in \bb{\mathbf{c}}
\}.
\]

On the other hand, since $(\mathfrak{P}^{(2)},s)<_{\mathbf{Pth}_{\boldsymbol{\mathcal{A}}^{(2)}}} (\mathfrak{Q}^{(2)},t)$ and the sequence witnessing this fact is entirely composed of coherent head-constant echelonless second-order paths we have, by Lemma~\ref{LDOrdV}, that 
\[
\max\{
\bb{T^{(1)}_{i}}\mid i\in \bb{\mathbf{c}}
\}
<
\max\{
\bb{U^{(1)}_{j}}\mid j\in \bb{\mathbf{d}}
\}
.
\]

This leads to the following contradiction
$$
\max\{
\bb{T^{(1)}_{i}}\mid i\in \bb{\mathbf{c}}
\}
<
\max\{
\bb{U^{(1)}_{j}}\mid j\in \bb{\mathbf{d}}
\}
<
\max\{
\bb{T^{(1)}_{i}}\mid i\in \bb{\mathbf{c}}
\}.
$$

Therefore, we can affirm that $\mathfrak{P}^{(2)}=\mathfrak{Q}^{(2)}$.

From this it follows that $\leq_{\coprod\mathrm{Pth_{\boldsymbol{\mathcal{A}}^{(2)}}}}$ is a partial order.

We next prove that no strictly decreasing $\omega_{0}$-chain can exist in the partially ordered set $(\coprod\mathrm{Pth}_{\boldsymbol{\mathcal{A}}^{(2)}}, \leq_{\mathbf{Pth}_{\boldsymbol{\mathcal{A}}^{(2)}}})$. Assume, towards a contradiction, that there exists an strictly decreasing $\omega_{0}$-chain $((\mathfrak{R}^{(2)}_{k},s_{k}))_{k\in\omega_{0}}$ in $(\coprod\mathrm{Pth}_{\boldsymbol{\mathcal{A}}^{(2)}},\leq_{\mathbf{Pth}_{\boldsymbol{\mathcal{A}}^{(2)}}})$. Then, for every $k\in\omega_{0}$, we have that 
\begin{enumerate}
\item $\mathfrak{R}^{(2)}_{k}$ is a second-order path in $\mathrm{Pth}_{\boldsymbol{\mathcal{A}}^{(2)},s_{k}}$, and
\item $(\mathfrak{R}^{(2)}_{k+1},s_{k+1})\prec_{\mathbf{Pth}_{\boldsymbol{\mathcal{A}}^{(2)}}}(\mathfrak{R}^{(2)}_{k},s_{k})$.
\end{enumerate}

Note that an strictly decreasing $\omega_{0}$-chain in $(\coprod\mathrm{Pth}_{\boldsymbol{\mathcal{A}}^{(2)}},\leq_{\mathbf{Pth}_{\boldsymbol{\mathcal{A}}^{(2)}}})$, in principle, should be defined by using, in clause (2) above, $<_{\mathbf{Pth}_{\boldsymbol{\mathcal{A}}^{(2)}}}$ instead of $\prec_{\mathbf{Pth}_{\boldsymbol{\mathcal{A}}^{(2)}}}$. However, by Remark~\ref{RDOrd}, such a chain would, ultimately, lead to a chain as above, where $\prec_{\mathbf{Pth}_{\boldsymbol{\mathcal{A}}^{(2)}}}$ is actually used.

By construction, for every $k\in\omega_{0}$, the labelled second-order path $(\mathfrak{R}^{(2)}_{k},s_{k})$ cannot be minimal in the partially ordered set $(\coprod\mathrm{Pth}_{\boldsymbol{\mathcal{A}}^{(2)}},\leq_{\mathbf{Pth}_{\boldsymbol{\mathcal{A}}^{(2)}}})$. Otherwise, by the minimality of $(\mathfrak{R}^{(2)}_{k},s_{k})$, the labelled second-order path $(\mathfrak{R}^{(2)}_{k+1},s_{k+1})$ would be equal to $(\mathfrak{R}^{(2)}_{k},s_{k})$, contradicting the fact that $\prec_{\mathbf{Pth}_{\boldsymbol{\mathcal{A}}^{(2)}}}$ and $\Delta_{\coprod\mathrm{Pth}_{\boldsymbol{\mathcal{A}}^{(2)}}}$ are disjoint relations. Thus, by Lemma~\ref{LDOrdI}, for every $k\in\omega_{0}$, we have that
\begin{enumerate}
\item $\mathfrak{R}^{(2)}_{k}$ is a $(2,[1])$-identity second-order path on a non-minimal path term class; or
\item $\mathfrak{R}^{(2)}_{k}$ is a second-order path of length strictly greater than one containing at least one   second-order echelon; or
\item $\mathfrak{R}^{(2)}_{k}$ is an echelonless second-order path.
\end{enumerate}

We claim that in the sequence $((\mathfrak{R}^{(2)}_{k},s_{k}))_{k\in\omega_{0}}$, no $(2,[1])$-identity second-order path can occur. Let us suppose, towards a contradiction, that this is not the case. Thus, assume that, for some $k\in\omega_{0}$, $\mathfrak{R}^{(2)}_{k}$ is a $(2,[1])$-identity second-order path, then the subsequence $((\mathfrak{R}^{(2)}_{l},s_{l}))_{l\geq k}$ is composed entirely of $(2,[1])$-identity second-order paths according to Proposition~\ref{PDULower}. Hence, for every $l\in\omega_{0}$ with $l\geq k$, the second-order path $\mathfrak{R}^{(2)}_{l}$ has the form
$$
\mathfrak{R}^{(2)}_{l}=\mathrm{ip}^{(2,[1])\sharp}_{s_{l}}\left(
[R_{l}]_{s_{l}}
\right),
$$
for a suitable path term class $[R_{l}]_{s_{l}}$ in $[\mathrm{PT}_{\boldsymbol{\mathcal{A}}}]_{s_{l}}$. 

Taking into account that we are assuming that, for every $l\in\omega_{0}$ with $l\geq k$,
$$
\left(\mathfrak{R}^{(2)}_{l+1},s_{l+1}
\right)
\prec_{\mathbf{Pth}_{\boldsymbol{\mathcal{A}}^{(2)}}}
\left(\mathfrak{R}^{(2)}_{l},s_{l}
\right)
$$
we have, according to Definition~\ref{DDOrd} that, for every $l\in\omega_{0}$ with $l\geq k$, the following inequality holds
$$
\left([R_{l+1}]_{s_{l+1}},s_{l+1}
\right)
<_{[\mathbf{PT}_{\boldsymbol{\mathcal{A}}}]}
\left([R_{l}]_{s_{l}},s_{l}
\right).
$$
Thus $(([R_{l}]_{s_{l}},s_{l}))_{l\geq k}$ is an strictly decreasing $\omega_{0}$-chain in 
$
(\coprod[\mathrm{PT}_{\boldsymbol{\mathcal{A}}}],
\leq_{[\mathbf{PT}_{\boldsymbol{\mathcal{A}}}]})
$, 
contradicting the fact that, by Proposition~\ref{PPTQOrdArt}, $
(\coprod[\mathrm{PT}_{\boldsymbol{\mathcal{A}}}],
\leq_{[\mathbf{PT}_{\boldsymbol{\mathcal{A}}}]})
$ is an Artinian ordered set.

Thus, for every $k\in\omega_{0}$, $\mathfrak{R}^{(2)}_{k}$ cannot be a $(2,[1])$-identity second-order path.

We claim that, in the sequence $((\mathfrak{R}^{(2)}_{k},s_{k}))_{k\in\omega_{0}}$, the number of second-order paths of length strictly greater than one containing at least one  -second-order echelon is bounded. Let us suppose towards a contradiction that this is not the case. Then we can extract from $((\mathfrak{R}^{(2)}_{k},s_{k}))_{k\in\omega_{0}}$ an infinite subsequence $((\mathfrak{R}^{(2)}_{\varphi(k)},s_{\varphi(k)}))_{k\in\omega_{0}}$, where $\varphi$ is an strictly increasing endomapping of $\omega_{0}$, such that, for every $k\in\omega_{0}$, we have that
\begin{enumerate}
\item $\mathfrak{R}^{(2)}_{\varphi(k)}$ is a second-order path of length strictly greater than one containing at least one   second-order echelon; and
\item $(\mathfrak{R}^{(2)}_{\varphi(k+1)}, s_{\varphi(k+1)})
<_{\mathbf{Pth}_{\boldsymbol{\mathcal{A}}^{(2)}}}
(\mathfrak{R}^{(2)}_{\varphi(k)}, s_{\varphi(k)})
$.
\end{enumerate}

But, by Corollary~\ref{CDOrdI}, this leads to a contradiction, concretely, to that of having an strictly decreasing $\omega_{0}$-chain in $(\mathbb{N},\leq)$ constructed from the length of such second-order paths.

From this, it follows that, in the sequence $((\mathfrak{R}^{(2)}_{k},s_{k}))_{k\in\omega_{0}}$, the number of indexes $k\in\omega_{0}$ for which $\mathfrak{R}^{(2)}_{k}$ is a second-order path of length strictly greater than one containing at least one   second-order echelon is bounded.  By a similar argument, the number of indexes $k\in\omega_{0}$ for which $\mathfrak{R}^{(2)}_{k}$ is an echelonless second-order path that is not head-constant is bounded. The same applies for the number of indexes $k\in\omega_{0}$ for which  $\mathfrak{R}^{(2)}_{k}$ is a head-constant echelonless second-order path that is not coherent.

Let $l\in\omega_{0}$ be the maximum of the indexes for which $\mathfrak{R}^{(2)}_{l}$ is either (1) a second-order path of length strictly greater than one containing at least one   second-order echelon, or (2) an echelonless second-order path that is not head-constant, or (3) head-constant echelonless second-order path that is not coherent. Then, for the infinite subsequence $((\mathfrak{R}^{(2)}_{k},s_{k}))_{k\in\omega_{0}-(l+1)}$ of $((\mathfrak{R}^{(2)}_{k},s_{k}))_{k\in\omega_{0}}$ and for every $k\in\omega_{0}-(l+1)$, we have that 
\begin{enumerate}
\item $\mathfrak{R}^{(2)}_{k}$ is a coherent head-constant echelonless second-order path; and
\item $(\mathfrak{R}^{(2)}_{k+1}, s_{k+1})
\prec_{\mathbf{Pth}_{\boldsymbol{\mathcal{A}}^{(2)}}}
(\mathfrak{R}^{(2)}_{k}, s_{k})
$.
\end{enumerate}
But, by Lemma~\ref{LDOrdV}, this leads to a contradiction, concretely to that of having an strictly decreasing $\omega_{0}$-chain in $(\mathbb{N},\leq)$ constructed from the maximum of the heights of the sequence of first-order  translations of such second-order paths.

Thus, no strictly decreasing $\omega_{0}$-chain can exist in $(\coprod\mathrm{Pth}_{\boldsymbol{\mathcal{A}}^{(2)}}, \leq_{\mathbf{Pth}_{\boldsymbol{\mathcal{A}}^{(2)}}})$. That is, the order $\leq_{\mathbf{Pth}_{\boldsymbol{\mathcal{A}}^{(2)}}}$ on $\coprod\mathrm{Pth}_{\boldsymbol{\mathcal{A}}^{(2)}}$ is Artinian.
\end{proof}

Taking into account the definitions above, now the coproduct of the $(2,[1])$-identity second-order path mapping becomes an order embedding from the Artinian ordered set $(\coprod[\mathrm{PT}_{\boldsymbol{\mathcal{A}}}], \leq_{[\mathbf{PT}_{\boldsymbol{\mathcal{A}}}]})$, introduced on Definition~\ref{DPTQOrd}, to the Artinian orderedset $(\coprod\mathrm{Pth}_{\boldsymbol{\mathcal{A}}^{(2)}}, \leq_{\mathbf{Pth}_{\boldsymbol{\mathcal{A}}^{(2)}}})$, introduced on Definition~\ref{DDOrd}.

\begin{restatable}{proposition}{PDUOrdEmb}
\label{PDUOrdEmb} The mapping $\coprod\mathrm{ip}^{(2,[1])\sharp}$ is an order embedding
$$
\textstyle
\coprod\mathrm{ip}^{(2,[1])\sharp}\colon
\left(
\coprod
[\mathrm{PT}_{\boldsymbol{\mathcal{A}}}], 
\leq_{[\mathbf{PT}_{\boldsymbol{\mathcal{A}}}]}
\right)
\mor 
\left(
\coprod\mathrm{Pth}_{\boldsymbol{\mathcal{A}}^{(2)}}, 
\leq_{\mathbf{Pth}_{\boldsymbol{\mathcal{A}}^{(2)}}}
\right).
$$
\end{restatable}
\begin{proof}
Let $s,t$ be sorts in $S$ and let $([Q]_{t},t)$ and let $([P]_{s},s)$ be pairs in $\coprod
[\mathrm{PT}_{\boldsymbol{\mathcal{A}}}]$. We need to prove that the following statements are equivalent
\begin{enumerate}
\item $([Q]_{t},t)\leq_{[\mathbf{PT}_{\boldsymbol{\mathcal{A}}}]} ([P]_{s},s)$;
\item $(\mathrm{ip}^{(2,[1])\sharp}_{t}([Q]_{t}),t)
\leq_{\mathbf{Pth}_{\boldsymbol{\mathcal{A}}^{(2)}}}
(\mathrm{ip}^{(2,[1])\sharp}_{s}([P]_{s}),s).
$ 
\end{enumerate}
Note that the statement trivially holds when the above statement are equalities. Therefore, it suffices to prove the above equivalence for the case of strict inequalities.

On the one hand, if  $([Q]_{t},t)<_{[\mathbf{PT}_{\boldsymbol{\mathcal{A}}}]} ([P]_{s},s)$, then, in virtue of Definition~\ref{DDOrd}, we have that   $(\mathrm{ip}^{(2,[1])\sharp}_{t}([Q]_{t}),t)
\prec_{\mathbf{Pth}_{\boldsymbol{\mathcal{A}}^{(2)}}}
(\mathrm{ip}^{(2,[1])\sharp}_{s}([P]_{s}),s).
$ 

On the other hand, if  $(\mathrm{ip}^{(2,[1])\sharp}_{t}([Q]_{t}),t)
<_{\mathbf{Pth}_{\boldsymbol{\mathcal{A}}^{(2)}}}
(\mathrm{ip}^{(2,[1])\sharp}_{s}([P]_{s}),s).
$ then, by Corollary~\ref{CDULower}, we have that   $([Q]_{t},t)<_{[\mathbf{PT}_{\boldsymbol{\mathcal{A}}}]} ([P]_{s},s)$.

This finishes the proof.
\end{proof}

Taking into account the definitions above, we next show that the coproduct of the $(2,0)$-identity second-order path mapping becomes an order embedding from the Artinian ordered set $(\coprod\mathrm{T}_{\Sigma}(X), \leq_{\mathbf{T}_{\Sigma}(X)})$, introduced on Remark~\ref{RTermOrd}, to the Artinian ordered set $(\coprod\mathrm{Pth}_{\boldsymbol{\mathcal{A}}^{(2)}}, \leq_{\mathbf{Pth}_{\boldsymbol{\mathcal{A}}^{(2)}}})$, introduced on Definition~\ref{DDOrd}.

\begin{restatable}{proposition}{PDZOrdEmb}
\label{PDZOrdEmb} The mapping $\coprod\mathrm{ip}^{(2,0)\sharp}$ is an order embedding
$$
\textstyle
\coprod\mathrm{ip}^{(2,0)\sharp}\colon
\left(
\coprod\mathrm{T}_{\Sigma}(X)
\leq_{\mathbf{T}_{\Sigma}(X)}
\right)
\mor 
\left(\coprod\mathrm{Pth}_{\boldsymbol{\mathcal{A}}^{(2)}}, \leq_{\mathbf{Pth}_{\boldsymbol{\mathcal{A}}^{(2)}}}
\right).
$$
\end{restatable}
\begin{proof}
Let $s,t$ be sorts in $S$ and let $(Q,t)$ and let $(P,s)$ be pairs in $\coprod\mathrm{T}_{\Sigma}(X)$. We need to prove that the following statements are equivalent
\begin{enumerate}
\item $(Q,t)\leq_{\mathbf{T}_{\Sigma}(X)} (P,s)$;
\item $(\mathrm{ip}^{(2,0)\sharp}_{t}(Q),t)
\leq_{\mathbf{Pth}_{\boldsymbol{\mathcal{A}}^{(2)}}}
(\mathrm{ip}^{(2,[1])\sharp}_{s}(P),s).
$ 
\end{enumerate}
Note that the following equivalences holds
\begin{flushleft}
$\left(Q,t
\right)\leq_{\mathbf{T}_{\Sigma}(X)} 
\left(P,s
\right)$
\allowdisplaybreaks
\begin{align*}
&\Longleftrightarrow
\left(\mathrm{ip}^{([1],0)\sharp}_{t}\left(Q
\right), t\right)
\leq_{[\mathbf{Pth}_{\boldsymbol{\mathcal{A}}}]}
\left(\mathrm{ip}^{([1],0)\sharp}_{s}\left(P
\right), s\right)
\tag{1}
\\&\Longleftrightarrow
\left(
\mathrm{CH}^{[1]}_{t}
\left(
\mathrm{ip}^{([1],0)\sharp}_{t}
\left(
Q
\right)\right), t
\right)
\leq_{[\mathbf{PT}_{\boldsymbol{\mathcal{A}}}]}
\left(
\mathrm{CH}^{[1]}_{s}
\left(
\mathrm{ip}^{([1],0)\sharp}_{s}
\left(
P
\right)\right), s
\right)
\tag{2}
\\&\Longleftrightarrow
\left(
\mathrm{ip}^{(2,[1])\sharp}_{t}
\left(
\mathrm{CH}^{[1]}_{t}
\left(
\mathrm{ip}^{([1],0)\sharp}_{t}
\left(Q
\right)\right)\right), t\right)
\leq_{[\mathbf{PT}_{\boldsymbol{\mathcal{A}}}]}
\\&\qquad\qquad\qquad\qquad\qquad\qquad\qquad\qquad
\left(
\mathrm{ip}^{(2,[1])\sharp}_{s}
\left(
\mathrm{CH}^{[1]}_{s}
\left(
\mathrm{ip}^{([1],0)\sharp}_{s}
\left(
P
\right)\right)\right), s
\right)
\tag{3}
\\&\Longleftrightarrow
\left(
\mathrm{ip}^{(2,0)\sharp}_{t}
\left(
Q\right),t\right)
\leq_{\mathbf{Pth}_{\boldsymbol{\mathcal{A}}^{(2)}}}
\left(
\mathrm{ip}^{(2,0)\sharp}_{s}
\left(P
\right),s\right).
\tag{4}
\end{align*}
\end{flushleft}

The first equivalence follows from the fact that, by Proposition~\ref{PCHOrdIp}, the mapping $\coprod\mathrm{ip}^{([1],0)\sharp}$ is an order embedding; the second equivalence follows from the fact that, by Theorem~\ref{TIsoOrd}, the mapping $\coprod\mathrm{CH}^{[1]}$ is an order isomorphism; the third equivalence follows from the fact that, by Proposition~\ref{PDUOrdEmb}, $\coprod\mathrm{ip}^{(2,[1])\sharp}$ is an order embedding; finally, the last equivalence follows from the fact that, according to Definition~\ref{DDScTgZ}, $\mathrm{ip}^{(2,0)\sharp}=
\mathrm{ip}^{(2,[1])\sharp}\circ\mathrm{CH}^{[1]}\circ\mathrm{ip}^{([1],0)\sharp}$.

This finishes the proof.
\end{proof}

\chapter{
\texorpdfstring
{The categorial signature $\Sigma^{\boldsymbol{\mathcal{A}}^{(2)}}$ determined by $\boldsymbol{\mathcal{A}}^{(2)}$}
{The categorial second-order signature}
}\label{S2D}

In this chapter, we introduce the categorical signature determined by $\boldsymbol{\mathcal{A}}^{(2)}$, denoted by $\Sigma^{\boldsymbol{\mathcal{A}}^{(2)}}$. This new signature extends the categorial signature $\Sigma^{\boldsymbol{\mathcal{A}}}$ by adding the operations of $1$-composition, $1$-source and $1$-target and, for every sort $s\in S$ and every second-order rewrite rule $\mathfrak{p}^{(2)}\in \mathcal{A}^{(2)}_{s}$, a new constant operation symbol associated to it. This extension allows us to compare the respective free algebras $\mathbf{T}_{\Sigma}(X)$, $\mathbf{T}_{\Sigma^{\boldsymbol{\mathcal{A}}}}(X)$ and $\mathbf{T}_{\Sigma^{\boldsymbol{\mathcal{A}}^{(2)}}}(X)$.
 It is shown the existence of an $\Sigma^{\boldsymbol{\mathcal{A}}}$-embedding $\eta^{(2,1)\sharp}$ from  $\mathbf{PT}_{\boldsymbol{\mathcal{A}}}$ to $\mathbf{T}_{\Sigma^{\boldsymbol{\mathcal{A}}^{(2)}}}(X)$, which makes it possible to view all the path terms  in $\mathbf{PT}_{\boldsymbol{\mathcal{A}}}$ as terms in $\mathbf{T}_{\Sigma^{\boldsymbol{\mathcal{A}}^{(2)}}}(X)$. Moreover, there exists a  $\Sigma$-embedding $\eta^{(2,0)\sharp}$ from  $\mathbf{T}_{\Sigma}(X)$ to $\mathbf{T}_{\Sigma^{\boldsymbol{\mathcal{A}}^{(2)}}}(X)$, which makes it possible to view all the terms  in $\mathbf{T}_{\Sigma}(X)$ as terms in $\mathbf{T}_{\Sigma^{\boldsymbol{\mathcal{A}}^{(2)}}}(X)$. The new categorical signature allows the consideration of a many-sorted partial $\Sigma^{\boldsymbol{\mathcal{A}}^{(2)}}$-algebraic structure on the set of many-sorted second-order paths, denoted by $\mathbf{Pth}_{\boldsymbol{\mathcal{A}}^{(2)}}$, with the natural interpretations of the new operation symbols. This chapter concludes with the definition, by Artinian recursion, of the second-order Curry-Howard mapping, denoted by $\mathrm{CH}^{(2)}$, which allows to assign, for each sort $s\in S$ and each second-order path $\mathfrak{P}^{(2)}$ in $\mathrm{Pth}_{\boldsymbol{\mathcal{A}}^{(2)},s}$, a term $\mathrm{CH}^{(2)}_{s}(\mathfrak{P}^{(2)})$ of sort $s$ in $\mathbf{T}_{\Sigma^{\boldsymbol{\mathcal{A}}^{(2)}}}(X)_{s}$. It is then shown that $\mathrm{CH}^{(2)}$ is a $\Sigma$-homomorphism but not a $\Sigma^{\boldsymbol{\mathcal{A}}}$-homomorphism nor a $\Sigma^{\boldsymbol{\mathcal{A}}^{(2)}}$-homomorphism. We also prove that $\mathrm{CH}^{(2)}$ is order-preserving.


We introduce the categorial signature determined  by $\mathcal{A}^{(2)}$ on $\Sigma$.

\begin{restatable}{definition}{DDSigmaCat}
\label{DDSigmaCat}
Let $\Sigma$ be an $S$-sorted signature. Then the \emph{categorial signature determined by $\mathcal{A}^{(2)}$ on $\Sigma$} 
\index{categorial signature!determined by $\mathcal{A}^{(2)}$, $\Sigma^{\boldsymbol{\mathcal{A}}^{(2)}}$}, denoted by $\Sigma^{\boldsymbol{\mathcal{A}}^{(2)}}$, is the $S$-sorted signature defined, for every $(\mathbf{s},s)\in S^{\star}\times S$, as follows:
$$
{\Sigma}^{\boldsymbol{\mathcal{A}}^{(2)}}_{\mathbf{s},s}=
\begin{cases}
{\Sigma}^{\boldsymbol{\mathcal{A}}}_{\mathbf{s},s}, 	
& \text{if $\mathbf{s}\not\in \{\lambda, (s),(s,s)\}$;}\\
{\Sigma}^{\boldsymbol{\mathcal{A}}}_{\mathbf{s},s}\amalg\mathcal{A}^{(2)}_{s}, 	
& \text{if $\mathbf{s}= \lambda$;}\\
{\Sigma}^{\boldsymbol{\mathcal{A}}}_{\mathbf{s},s}\amalg\{\mathrm{sc}^{1}_{s}, \mathrm{tg}^{1}_{s}\}, 	
& \text{if $\mathbf{s}= (s)$;}\\
{\Sigma}^{\boldsymbol{\mathcal{A}}}_{\mathbf{s},s}\amalg\{\circ^{1}_{s}\}, 	
& \text{if $\mathbf{s}= (s,s)$.}
\end{cases}
$$
That is, ${\Sigma}^{\boldsymbol{\mathcal{A}}^{(2)}}$ is the expansion of $\Sigma^{\boldsymbol{\mathcal{A}}}$ obtained by adding, for every sort $s\in S$,
(1) as many constants of coarity $s$ as there are second-order rewrite rules in $\mathcal{A}^{(2)}_{s}$, (2) two unary operation symbols $\mathrm{sc}^{1}_{s}$ and $\mathrm{tg}^{1}_{s}$, both of them of arity $(s)$ and coarity $s$, which will be interpreted as total unary operations, and (3) a binary operation symbol $\circ^{1}_{s}$, of arity $(s,s)$ and coarity $s$, which will be interpreted as a binary (partial or total, depending on the case at hand) operation. 

Let us recall from Definition~\ref{DSigmaCat}, that ${\Sigma}^{\boldsymbol{\mathcal{A}}}$ was the expansion of $\Sigma$ obtained by adding, for every sort $s\in S$,
(1) as many constants of coarity $s$ as there are rewrite rules in $\mathcal{A}_{s}$, (2) two unary operation symbols $\mathrm{sc}^{0}_{s}$ and $\mathrm{tg}^{0}_{s}$, both of them of arity $(s)$ and coarity $s$, which will be interpreted as total unary operations, and (3) a binary operation symbol $\circ^{0}_{s}$, of arity $(s,s)$ and coarity $s$, which will be interpreted as a binary (partial or total, depending on the case at hand) operation.
\end{restatable}

We now define several $S$-sorted mappings from (1) the $S$-sorted set of variables $X$, (2) the underlying $S$-sorted set of rewrite rules $\mathcal{A}$ of $\boldsymbol{\mathcal{A}}$ and (3) the underlying $S$-sorted set of second-order rewrite rules $\mathcal{A}^{(2)}$ of $\boldsymbol{\mathcal{A}}^{(2)}$, into the underlying $S$-sorted set of the free $\Sigma^{\boldsymbol{\mathcal{A}}^{(2)}}$-algebra on $X$ (we recall that $X$ underlies $\boldsymbol{\mathcal{A}}$ and $\boldsymbol{\mathcal{A}}^{(2)}$ and that the latter have already been fixed).  

\begin{restatable}{definition}{DDEta}
\label{DDEta}
\index{terms!second-order!$\mathbf{T}_{\Sigma^{\boldsymbol{\mathcal{A}}^{(2)}}}(X)$}
Let $\mathbf{T}_{\Sigma^{\boldsymbol{\mathcal{A}}^{(2)}}}(X)$ be the free $\Sigma^{\boldsymbol{\mathcal{A}}^{(2)}}$-algebra on $X$. We will denote by
\begin{enumerate}
\item $\eta^{(2,X)}$ the $S$-sorted mapping from $X$ to 
$\mathrm{T}_{\Sigma^{\boldsymbol{\mathcal{A}}^{(2)}}}(X)$ such that, for every sort $s\in S$, sends an element $x\in X_{s}$ to the variable  $x\in \mathrm{T}_{\Sigma^{\boldsymbol{\mathcal{A}}^{(2)}}}(X)_{s}$; by
\index{inclusion!second-order!$\eta^{(2,X)}$}
\item $\eta^{(2,\mathcal{A})}$ the $S$-sorted mapping from 
$\mathcal{A}$ to 
$\mathrm{T}_{\Sigma^{\boldsymbol{\mathcal{A}}^{(2)}}}(X)$ 
such that, for every sort $s\in S$, sends a  rewrite rule $\mathfrak{p}\in \mathcal{A}_{s}$ to the constant $\mathfrak{p}^{\mathrm{T}_{\Sigma^{\boldsymbol{\mathcal{A}}^{(2)}}}(X)}$, and by
\index{inclusion!second-order!$\eta^{(2,\mathcal{A})}$}
\item $\eta^{(2,\mathcal{A}^{(2)})}$ the $S$-sorted mapping from 
$\mathcal{A}^{(2)}$ to 
$\mathrm{T}_{\Sigma^{\boldsymbol{\mathcal{A}}^{(2)}}}(X)$ 
such that, for every sort $s\in S$, sends a second-order rewrite rule $\mathfrak{p}^{(2)}\in \mathcal{A}^{(2)}_{s}$ to the constant $\mathfrak{p}^{(2)\mathrm{T}_{\Sigma^{\boldsymbol{\mathcal{A}}^{(2)}}}(X)}$.
\index{inclusion!second-order!$\eta^{(2,\mathcal{A}^{(2)})}$}
\end{enumerate}
The above $S$-sorted mappings are depicted in the diagram of Figure~\ref{FDPthEmb}.
\end{restatable}

\begin{figure}
\begin{tikzpicture}
[ACliment/.style={-{To [angle'=45, length=5.75pt, width=4pt, round]}},scale=.8]
\node[] (x) at (0,0) [] {$X$};
\node[] (a) at (0,-1.5) [] {$\mathcal{A}$};
\node[] (a2) at (0,-3) [] {$\mathcal{A}^{(2)}$};
\node[] (T) at (6,-3) [] {$\mathrm{T}_{\Sigma^{\boldsymbol{\mathcal{A}}^{(2)}}}(X)$};
\draw[ACliment, bend left=20]  (x) to node [above right] {$\eta^{(2,X)}$} (T);
\draw[ACliment, bend left=10]  (a) to node [midway, fill=white] {$\eta^{(2,\mathcal{A})}$} (T);
\draw[ACliment]  (a2) to node [below] {$\eta^{(2,\mathcal{A}^{(2)})}$} (T);
\end{tikzpicture}
\caption{Embeddings relative to $X$, $\mathcal{A}$ and $\mathcal{A}^{(2)}$ at layer 2.}\label{FDPthEmb}
\end{figure}

We next investigate how the just defined embeddings interact with the free algebras of previous layers. In this regard, we next define the $\Sigma^{\boldsymbol{\mathcal{A}}}$-reduct of the free $\Sigma^{\boldsymbol{\mathcal{A}}^{(2)}}$-algebra $\mathbf{T}_{\Sigma^{\boldsymbol{\mathcal{A}}^{(2)}}}(X)$.

\begin{definition}\label{DDURed}
Let $\mathrm{in}^{\Sigma,(2,1)}$ be the canonical embedding of $\Sigma^{\boldsymbol{\mathcal{A}}}$ into
$\Sigma^{\boldsymbol{\mathcal{A}}^{(2)}}$. Then, by Proposition~\ref{PFunSig}, for the morphism $\mathbf{in}^{\Sigma,(2,1)} = (\mathrm{id}^{S},\mathrm{in}^{\Sigma,(2,1)})$ from $(S,\Sigma^{\boldsymbol{\mathcal{A}}})$ to
$(S,\Sigma^{\boldsymbol{\mathcal{A}}^{(2)}})$ and the free $\Sigma^{\boldsymbol{\mathcal{A}}^{(2)}}$-algebra $\mathbf{T}_{\Sigma^{\boldsymbol{\mathcal{A}}^{(2)}}}(X)$, we will denote by $\mathbf{T}_{\Sigma^{\boldsymbol{\mathcal{A}}^{(2)}}}^{(1,2)}(X)$ the $\Sigma^{\boldsymbol{\mathcal{A}}}$-algebra $(\mathbf{in}^{\Sigma,(1,2)})(\mathbf{T}_{\Sigma^{\boldsymbol{\mathcal{A}}^{(2)}}}(X))$. We will call $\mathbf{T}_{\Sigma^{\boldsymbol{\mathcal{A}}^{(2)}}}^{(1,2)}(X)$ the \emph{$\Sigma^{\boldsymbol{\mathcal{A}}}$-reduct} of the free $\Sigma^{\boldsymbol{\mathcal{A}}^{(2)}}$-algebra $\mathbf{T}_{\Sigma^{\boldsymbol{\mathcal{A}}^{(2)}}}(X)$.
\end{definition}

\begin{remark}\label{RDURed} The underlying $S$-sorted set of $\mathbf{T}_{\Sigma^{\boldsymbol{\mathcal{A}}^{(2)}}}^{(1,2)}(X)$ is the same as that of the free $\Sigma^{\boldsymbol{\mathcal{A}}^{(2)}}$-algebra $\mathbf{T}_{\Sigma^{\boldsymbol{\mathcal{A}}^{(2)}}}(X)$, while, for every pair $(\mathbf{s},s)\in S^{\star}\times S$ and every operation symbol $\tau\in\Sigma^{\boldsymbol{\mathcal{A}}}_{\mathbf{s},s}$,
$$
\tau^{\mathbf{T}^{(1,2)}_{\Sigma^{\boldsymbol{\mathcal{A}}^{(2)}}}(X)}=\tau^{\mathbf{T}_{\Sigma^{\boldsymbol{\mathcal{A}}^{(2)}}}(X)}.
$$
\end{remark}

\begin{restatable}{proposition}{PDUEmb}
\label{PDUEmb}
\index{inclusion!second-order!$\eta^{(2,1)\sharp}$}
There exists an embedding from the free partial $\Sigma^{\boldsymbol{\mathcal{A}}}$-algebra $\mathbf{PT}_{\boldsymbol{\mathcal{A}}}$ into the $\Sigma^{\boldsymbol{\mathcal{A}}}$-reduct $\mathbf{T}^{(1,2)}_{\Sigma^{\boldsymbol{\mathcal{A}}^{(2)}}}(X)$.
\end{restatable}
\begin{proof}
Consider $\eta^{(1,X)}$ and $\eta^{(2,X)}$, the canonical embeddings of $X$ into $\mathrm{T}_{\Sigma^{\boldsymbol{\mathcal{A}}}}(X)$ and $\mathrm{T}_{\Sigma^{\boldsymbol{\mathcal{A}}^{(2)}}}(X)$, respectively. Since $\mathbf{T}^{(1,2)}_{\Sigma^{\boldsymbol{\mathcal{A}}^{(2)}}}(X)$ is a $\Sigma^{\boldsymbol{\mathcal{A}}}$-algebra, by the universal property of the free $\Sigma^{\boldsymbol{\mathcal{A}}}$-algebra on $X$, i.e., $\mathbf{T}_{\Sigma^{\boldsymbol{\mathcal{A}}}}(X)$, there exists a unique many-sorted $\Sigma^{\boldsymbol{\mathcal{A}}}$-homomorphism, that we will denote by $\eta^{(2,1)\sharp}$ from $\mathbf{T}_{\Sigma^{\boldsymbol{\mathcal{A}}}}(X)$ to $\mathbf{T}^{(1,2)}_{\Sigma^{\boldsymbol{\mathcal{A}}^{(2)}}}(X)$ such that 
$$
\eta^{(2,1)\sharp}\circ\eta^{(1,X)}=\eta^{(2,X)}.
$$

The many-sorted $\Sigma^{\boldsymbol{\mathcal{A}}}$-homomorphism $\eta^{(2,1)\sharp}$ sends, for every sort $s\in S$, a term $P\in\mathrm{T}_{\Sigma^{\boldsymbol{\mathcal{A}}}}(X)_{s}$ to the term $P\in \mathrm{T}_{\Sigma^{\boldsymbol{\mathcal{A}}^{(2)}}}(X)_{s}$. This many-sorted mapping is injective.

At this point, we can, respectively, correstrict $\eta^{(1,X)}$ and restrict $\eta^{(2,1)\sharp}$, to the many-sorted set of path terms, i.e., to $\mathrm{PT}_{\boldsymbol{\mathcal{A}}}$. To obtain the desired result. To simplify the presentation, we will keep using $\eta^{(1,X)}$ and $\eta^{(2,1)\sharp}$. 

As usual, we identify $\mathrm{PT}_{\boldsymbol{\mathcal{A}}}$ with its image under $\eta^{(2,1)\sharp}$ in $\mathrm{T}_{\Sigma^{\boldsymbol{\mathcal{A}}^{(2)}}}(X)$. In this way, $\mathrm{PT}_{\boldsymbol{\mathcal{A}}}$ becomes a subset of $\mathrm{T}_{\Sigma^{\boldsymbol{\mathcal{A}}^{(2)}}}(X)$ and $\mathbf{PT}_{\boldsymbol{\mathcal{A}}}$ a subalgebra of $\mathbf{T}_{\Sigma^{\boldsymbol{\mathcal{A}}^{(2)}}}^{(1,2)}(X)$.
\end{proof}

Stepping down one more layer, we now introduce the $\Sigma$-reduct of the free $\Sigma^{\boldsymbol{\mathcal{A}}^{(2)}}$-algebra $\mathbf{T}_{\Sigma^{\boldsymbol{\mathcal{A}}^{(2)}}}(X)$.

\begin{definition}\label{DDZRed}
Let $\mathrm{in}^{\Sigma,(2,0)}$ be the canonical embedding of $\Sigma$ into
$\Sigma^{\boldsymbol{\mathcal{A}}^{(2)}}$. Then, by Proposition~\ref{PFunSig}, for the morphism $\mathbf{in}^{\Sigma,(2,0)} = (\mathrm{id}^{S},\mathrm{in}^{\Sigma,(2,0)})$ from $(S,\Sigma)$ to
$(S,\Sigma^{\boldsymbol{\mathcal{A}}^{(2)}})$ and the free $\Sigma^{\boldsymbol{\mathcal{A}}^{(2)}}$-algebra $\mathbf{T}_{\Sigma^{\boldsymbol{\mathcal{A}}^{(2)}}}(X)$, we will denote by $\mathbf{T}_{\Sigma^{\boldsymbol{\mathcal{A}}^{(2)}}}^{(0,2)}(X)$ the $\Sigma$-algebra $\mathbf{in}^{\Sigma,(0,2)}(\mathbf{T}_{\Sigma^{\boldsymbol{\mathcal{A}}^{(2)}}}(X))$. We will call $\mathbf{T}_{\Sigma^{\boldsymbol{\mathcal{A}}^{(2)}}}^{(0,2)}(X)$ the \emph{$\Sigma$-reduct} of the free $\Sigma^{\boldsymbol{\mathcal{A}}^{(2)}}$-algebra $\mathbf{T}_{\Sigma^{\boldsymbol{\mathcal{A}}^{(2)}}}(X)$.

Let us note that $\mathrm{in}^{\Sigma,(2,0)}$ factors through the mappings $\mathrm{in}^{\Sigma,(2,1)}$ and $\mathrm{in}^{\Sigma,(1,0)}$, i.e., the following equation holds
$$
\mathrm{in}^{\Sigma,(2,0)}=\mathrm{in}^{\Sigma,(2,1)}\circ\mathrm{in}^{\Sigma,(1,0)}.
$$
\end{definition}

\begin{remark}\label{RDZRed} The underlying $S$-sorted set of $\mathbf{T}_{\Sigma^{\boldsymbol{\mathcal{A}}^{(2)}}}^{(0,2)}(X)$ is the same as that of the free $\Sigma^{\boldsymbol{\mathcal{A}}^{(2)}}$-algebra $\mathbf{T}_{\Sigma^{\boldsymbol{\mathcal{A}}^{(2)}}}(X)$, while, for every pair $(\mathbf{s},s)\in S^{\star}\times S$ and every operation symbol $\sigma\in\Sigma_{\mathbf{s},s}$,
$$
\sigma^{\mathbf{T}^{(0,2)}_{\Sigma^{\boldsymbol{\mathcal{A}}^{(2)}}}(X)}=\sigma^{\mathbf{T}_{\Sigma^{\boldsymbol{\mathcal{A}}^{(2)}}}(X)}.
$$
\end{remark}

\begin{restatable}{proposition}{PDZEmb}
\label{PDZEmb}
\index{inclusion!second-order!$\eta^{(2,0)\sharp}$}
There exists an embedding from the free $\Sigma$-algebra $\mathbf{T}_{\Sigma}(X)$ into the $\Sigma$-reduct $\mathbf{T}^{(0,2)}_{\Sigma^{\boldsymbol{\mathcal{A}}^{(2)}}}(X)$.
\end{restatable}

\begin{proof}
Consider $\eta^{(0,X)}$ and $\eta^{(2,X)}$, the canonical embeddings of $X$ into $\mathrm{T}_{\Sigma}(X)$ and $\mathrm{T}_{\Sigma^{\boldsymbol{\mathcal{A}}^{(2)}}}(X)$, respectively. Since $\mathbf{T}^{(0,2)}_{\Sigma^{\boldsymbol{\mathcal{A}}^{(2)}}}(X)$ is a $\Sigma$-algebra, by the universal property of the free $\Sigma$-algebra on $X$, i.e., $\mathbf{T}_{\Sigma}(X)$, there exists a unique many-sorted $\Sigma$-homomorphism, that we will denote by $\eta^{(2,0)\sharp}$ from $\mathbf{T}_{\Sigma}(X)$ to $\mathbf{T}^{(0,2)}_{\Sigma^{\boldsymbol{\mathcal{A}}^{(2)}}}(X)$ such that 
$$
\eta^{(2,0)\sharp}\circ\eta^{(0,X)}=\eta^{(2,X)}.
$$

The many-sorted $\Sigma$-homomorphism $\eta^{(2,0)\sharp}$ sends, for every sort $s\in S$, a term $P\in\mathrm{T}_{\Sigma}(X)_{s}$ to the term $P\in \mathrm{T}_{\Sigma^{\boldsymbol{\mathcal{A}}^{(2)}}}(X)_{s}$. This many-sorted mapping is injective. The reader is advised to consult the diagram presented in Figure~\ref{FDAEmb}.

As usual, we identify $\mathrm{T}_{\Sigma}(X)$ with its image under $\eta^{(2,0)\sharp}$ in $\mathrm{T}_{\Sigma^{\boldsymbol{\mathcal{A}}^{(2)}}}(X)$. In this way, $\mathrm{T}_{\Sigma}(X)$ becomes a subset of $\mathrm{T}_{\Sigma^{\boldsymbol{\mathcal{A}}^{(2)}}}(X)$ and $\mathbf{T}_{\Sigma}(X)$ a subalgebra of $\mathbf{T}_{\Sigma^{\boldsymbol{\mathcal{A}}^{(2)}}}^{(0,2)}(X)$.
\end{proof}

\begin{proposition}\label{PDUZEmb}
There exists an embedding from the $\Sigma$-reduct $\mathbf{PT}^{(0,1)}_{\boldsymbol{\mathcal{A}}}$ into the $\Sigma$-reduct $\mathbf{T}^{(0,2)}_{\Sigma^{\boldsymbol{\mathcal{A}}^{(2)}}}(X)$.
\end{proposition}
\begin{proof}
The $\Sigma^{\boldsymbol{\mathcal{A}}}$-embedding $\eta^{(2,1)\sharp}$ from $\mathbf{PT}_{\boldsymbol{\mathcal{A}}}$ to $\mathbf{T}^{(1,2)}_{\Sigma^{\boldsymbol{\mathcal{A}}^{(2)}}}(X)$ presented in Proposition~\ref{PDUEmb} is, in particular, a $\Sigma$-embedding from $\mathbf{PT}^{(0,1)}_{\boldsymbol{\mathcal{A}}}$ to $\mathbf{T}^{(0,2)}_{\Sigma^{\boldsymbol{\mathcal{A}}^{(2)}}}(X)$.
\end{proof}

The embedding of the free $\Sigma$-algebra $\mathbf{T}_{\Sigma}(X)$ into the $\Sigma$-reduct $\mathbf{T}^{(0,2)}_{\Sigma^{\boldsymbol{\mathcal{A}}^{(2)}}}(X)$ factors through the embedding of  $\mathbf{T}_{\Sigma}(X)$ into the $\Sigma$-reduct $\mathbf{PT}^{(0,1)}_{\boldsymbol{\mathcal{A}}}$ presented in Proposition~\ref{PEmb} and the embedding of $\Sigma$-reduct $\mathbf{PT}^{(0,1)}_{\boldsymbol{\mathcal{A}}}$ into the $\Sigma$-reduct $\mathbf{T}^{(0,2)}_{\Sigma^{\boldsymbol{\mathcal{A}}^{(2)}}}(X)$ presented in Proposition~\ref{PDUZEmb}
.

\begin{proposition}\label{PDEmb} The equality $\eta^{(2,0)\sharp}=\eta^{(2,1)\sharp}\circ\eta^{(1,0)\sharp}$ holds.
\end{proposition}
\begin{proof}
This proposition entails that the diagram in Figure~\ref{FDEmb} commutes.
Note that $\eta^{(2,1)\sharp}\circ\eta^{(1,0)\sharp}$ is a $\Sigma$-homomorphism. Moreover, the following chain of equalities holds
\begin{align*}
\eta^{(2,1)\sharp}\circ\eta^{(1,0)\sharp}\circ\eta^{(0,X)}
&=
\eta^{(2,1)\sharp}\circ\eta^{(1,X)}
\tag{1}
\\&=
\eta^{(2,X)}.\tag{2}
\end{align*}
The first equation follows from the universal property of $\eta^{(1,0)\sharp}$ presented in Proposition~\ref{PEmb}; the last equation follows from the universal property of $\eta^{(2,1)\sharp
}$ presented in Proposition~\ref{PDUEmb}.

All in all, we can affirm that $\eta^{(2,0)\sharp}=\eta^{(2,1)\sharp}\circ\eta^{(1,0)\sharp}$.
\end{proof}

\begin{figure}
\begin{center}
\begin{tikzpicture}
[ACliment/.style={-{To [angle'=45, length=5.75pt, width=4pt, round]}},scale=0.8]
\node[] (x) at (0,0) [] {$X$};
\node[] (t) at (6,0) [] {$\mathrm{T}_{\Sigma}(X)$};
\node[] (t1) at (6,-2) [] {$\mathrm{PT}^{(0,1)}_{\boldsymbol{\mathcal{A}}}$};
\node[] (t2) at (6,-4) [] {$\mathrm{T}^{(0,2)}_{\Sigma^{\boldsymbol{\mathcal{A}}^{(2)}}}(X)$};
\node[] () at 	(9,-2) 	[] 	{$\textstyle \eta^{(2,0)\sharp}$};

\draw[ACliment]  (x) to node [ above right]
{$\textstyle \eta^{(0,X)}$} (t);
\draw[ACliment, bend right=10]  (x) to node [midway, fill=white]
{$\textstyle \eta^{(1,X)}$} (t1);
\draw[ACliment, bend right=20]  (x) to node [ below left]
{$\textstyle \eta^{(2,X)}$} (t2);

\draw[ACliment] 
($(t)+(.15,-.35)$) to node [right] {$\textstyle \eta^{(1,0)\sharp}$}  ($(t1)+(.15,.35)$);
\draw[ACliment] 
($(t1)+(.15,-.35)$) to node [right] {$\textstyle \eta^{(2,1)\sharp}$}  ($(t2)+(.15,.35)$);

\draw[ACliment, rounded corners] (t.east)
--
 ($(t.east)+(1.15,0)$)
--
($(t.east)+(1.15,-4)$)
-- (t2.east);
\end{tikzpicture}
\end{center}
\caption{Embeddings relative to $X$ at layers 0, 1 \& 2.}
\label{FDEmb}
\end{figure}

\begin{proposition}\label{PDAEmb} The equality $\eta^{(2,\mathcal{A})}=\eta^{(2,1)\sharp}\circ\eta^{(1,\mathcal{A})}$ holds.
\end{proposition}

\begin{figure}
\begin{center}
\begin{tikzpicture}
[ACliment/.style={-{To [angle'=45, length=5.75pt, width=4pt, round]}
}, scale=0.8]
\node[] (xq) 		at 	(0,0) 	[] 	{$\mathcal{A}$};
\node[] (txq) 	at 	(6,0) 	[] 	{$\mathrm{PT}_{\boldsymbol{\mathcal{A}}}$};
\node[] (txqc) 	at 	(6,-3) 	[] 	{$\mathrm{T}_{\Sigma^{\boldsymbol{\mathcal{A}}^{(2)}}}
(X)$};
\draw[ACliment]  (xq) 	to node [above right]	
{$\eta^{(1,\mathcal{A})}$} (txq);
\draw[ACliment]  (txq) 	to node [right]	
{$\eta^{(2,1)\sharp}$} (txqc);
\draw[ACliment, bend right=10]  (xq) 	to node [below left]	
{$\eta^{(2,\mathcal{A})}$} (txqc);
\end{tikzpicture}
\end{center}
\caption{Embeddings relative to $\mathcal{A}$ at layers 1 \& 2.}
\label{FDAEmb}
\end{figure}

\section{
\texorpdfstring
{A structure of partial $\Sigma^{\boldsymbol{\mathcal{A}}^{(2)}}$-algebra on $\mathrm{Pth}_{\boldsymbol{\mathcal{A}}^{(2)}}$}
{A partial algebra on second-order paths}
}
We next show that the $S$-sorted set $\mathrm{Pth}_{\boldsymbol{\mathcal{A}}^{(2)}}$ has a natural structure of many-sorted partial $\Sigma^{\boldsymbol{\mathcal{A}}^{(2)}}$-algebra.

\begin{restatable}{proposition}{PDPthDCatAlg}
\label{PDPthDCatAlg}
\index{path!second-order!$\mathbf{Pth}_{\boldsymbol{\mathcal{A}}^{(2)}}$}
The $S$-sorted set $\mathrm{Pth}_{\boldsymbol{\mathcal{A}}^{(2)}}$ is equipped, in a natural way, with a structure of many-sorted  partial $\Sigma^{\boldsymbol{\mathcal{A}}^{(2)}}$-algebra.
\end{restatable}

\begin{proof}
Let us denote by $\mathbf{Pth}_{\boldsymbol{\mathcal{A}}^{(2)}}$ the $\Sigma^{\boldsymbol{\mathcal{A}}^{(2)}}$-algebra defined as follows:

\textsf{(1)} The underlying $S$-sorted set of $\mathbf{Pth}_{\boldsymbol{\mathcal{A}}^{(2)}}$ is $\mathrm{Pth}_{\boldsymbol{\mathcal{A}}^{(2)}} = (\mathrm{Pth}_{\boldsymbol{\mathcal{A}}^{(2)},s})_{s\in S}$.

\textsf{(2)} For every $(\mathbf{s},s)\in S^{\star}\times S$ and every operation symbol $\sigma\in\Sigma_{\mathbf{s},s}$, the operation $\sigma^{\mathbf{Pth}_{\boldsymbol{\mathcal{A}}^{(2)}}}$ is given by the interpretation of $\sigma$ in $\mathbf{Pth}^{(0,2)}_{\boldsymbol{\mathcal{A}}^{(2)}}$ that, we recall, was stated in Proposition~\ref{PDPthAlg}, where, in addition, we proved that $\sigma^{\mathbf{Pth}_{\boldsymbol{\mathcal{A}}^{(2)}}}$ is well-defined. 

\textsf{(3)} For every $s\in S$ and every $\mathfrak{p}\in\mathcal{A}_{s}$, the constant $\mathfrak{p}^{\mathbf{Pth}_{\boldsymbol{\mathcal{A}}^{(2)}}}$ is given by the interpretation of $\mathfrak{p}$ in $\mathbf{Pth}^{(1,2)}_{\boldsymbol{\mathcal{A}}^{(2)}}$ that, we recall, was stated in Proposition~\ref{PDPthCatAlg}, where, in addition, we proved that $\mathfrak{p}^{\mathbf{Pth}_{\boldsymbol{\mathcal{A}}^{(2)}}}$ is well-defined. Let us recall that the constant $\mathfrak{p}^{\mathbf{Pth}_{\boldsymbol{\mathcal{A}}^{(2)}}}$ is given by
$$
\mathfrak{p}^{\mathbf{Pth}_{\boldsymbol{\mathcal{A}}^{(2)}}}
=
\mathrm{ech}^{(2,\mathcal{A})}_{s}\left(
\mathfrak{p}
\right)=
\mathrm{ip}^{(2,[1])\sharp}_{s}\left(
\eta^{([1],\mathcal{A})}_{s}\left(
\mathfrak{p}
\right)\right)
=
\left(\left[\mathfrak{p}^{\mathbf{PT}_{\boldsymbol{\mathcal{A}}}}
\right]_{s},\lambda,\lambda
\right).
$$

\textsf{(4)} For every $s\in S$, the unary operation of $0$-source, i.e., $\mathrm{sc}_{s}^{0\mathbf{Pth}_{\boldsymbol{\mathcal{A}}^{(2)}}}$, is given by the interpretation of $\mathrm{sc}_{s}^{0}$ in $\mathbf{Pth}^{(1,2)}_{\boldsymbol{\mathcal{A}}^{(2)}}$ that, we recall, was stated in Proposition~\ref{PDPthCatAlg}, where, in addition, we proved that $\mathrm{sc}_{s}^{0\mathbf{Pth}_{\boldsymbol{\mathcal{A}}^{(2)}}}$ is well-defined. 

\textsf{(5)} For every $s\in S$, the unary operation of $0$-target, i.e., $\mathrm{tg}_{s}^{0\mathbf{Pth}_{\boldsymbol{\mathcal{A}}^{(2)}}}$, is given by the interpretation of $\mathrm{tg}_{s}^{0}$ in $\mathbf{Pth}^{(1,2)}_{\boldsymbol{\mathcal{A}}^{(2)}}$ that, we recall, was stated in Proposition~\ref{PDPthCatAlg}, where, in addition, we proved that $\mathrm{tg}_{s}^{0\mathbf{Pth}_{\boldsymbol{\mathcal{A}}^{(2)}}}$ is well-defined. 

\textsf{(6)} For every $s\in S$, the partial binary operation of $0$-composition, i.e., $\circ_{s}^{0\mathbf{Pth}_{\boldsymbol{\mathcal{A}}^{(2)}}}$, is given by the interpretation of $\circ_{s}^{0}$ in $\mathbf{Pth}^{(1,2)}_{\boldsymbol{\mathcal{A}}^{(2)}}$ that, we recall, was stated in Proposition~\ref{PDPthCatAlg}, where, in addition, we proved that, given two second-order paths $\mathfrak{P}^{(2)}$ and $\mathfrak{Q}^{(2)}$ satisfying that 
$$
\mathrm{sc}^{(0,2)}_{s}
\left(
\mathfrak{Q}^{(2)}
\right)=
\mathrm{tg}^{(0,2)}_{s}
\left(\mathfrak{P}^{(2)}
\right),
$$
then the $0$-composition  $\mathfrak{Q}^{(2)}\circ_{s}^{0\mathbf{Pth}_{\boldsymbol{\mathcal{A}}^{(2)}}}\mathfrak{P}^{(2)}$ is well-defined.

\textsf{(7)} For every $s\in S$ and every $\mathfrak{p}^{(2)}\in\mathcal{A}^{(2)}_{s}$ with $\mathfrak{p}^{(2)}=([M]_{s},[N]_{s})$, the constant $\mathfrak{p}^{(2)\mathbf{Pth}_{\boldsymbol{\mathcal{A}}^{(2)}}}$ is given by the second-order echelon determined by $\mathfrak{p}^{(2)}$, i.e., by the one-step second-order path of sort $s$, i.e., $\mathrm{ech}^{(2,\mathcal{A})^{(2)}}_{s}(\mathfrak{p}^{(2)})$, that has $[M]_{s}$ as $([1],2)$-source and  $[N]_{s}$ as $([1],2)$-target. That is, the constant $\mathfrak{p}^{(2)\mathbf{Pth}_{\boldsymbol{\mathcal{A}}^{(2)}}}$ is interpreted as the following second-order echelon:
\begin{center}
\begin{tikzpicture}
[ACliment/.style={-{To [angle'=45, length=5.75pt, width=4pt, round]},font=\scriptsize},
AClimentD/.style={double equal sign distance, -implies, font=\scriptsize}
]
\node[] (1) at (0,0) [] {$\mathrm{ech}^{(2,\mathcal{A})^{(2)}}_{s}\left(
\mathfrak{p}^{(2)}
\right)\colon [M]_{s}$};
\node[] (2) at (6,0) [] {$[N]_{s}$.};
\draw[AClimentD]  (1) to node [above]
{$(\mathfrak{p}^{(2)},\mathrm{id}^{\mathrm{T}_{\Sigma^{\boldsymbol{\mathcal{A}}}}(X)_{s}})$} (2);
\end{tikzpicture}
\end{center}

\textsf{(8)} For every $s\in S$, the operation of $1$-source, i.e., $\mathrm{sc}_{s}^{1\mathbf{Pth}_{\boldsymbol{\mathcal{A}}^{(2)}}}$, is equal to 
$
\mathrm{ip}^{(2,[1])\sharp}_{s}\circ
\mathrm{sc}^{([1],2)}_{s}
$, i.e., to the mapping  that assigns to a second-order path $\mathfrak{P}^{(2)}$ in $\mathrm{Pth}_{\boldsymbol{\mathcal{A}}^{(2)}, s}$ the $(2,[1])$-identity  second-order path on the $([1],2)$-source of $\mathfrak{P}^{(2)}$.

\textsf{(9)} For every $s\in S$, the operation of $1$-target, i.e., $\mathrm{tg}_{s}^{1\mathbf{Pth}_{\boldsymbol{\mathcal{A}}^{(2)}}}$, is equal to 
$
\mathrm{ip}^{(2,[1])\sharp}_{s}\circ
\mathrm{tg}^{([1],2)}_{s}
$
, i.e., to the mapping  that assigns to a second-order path $\mathfrak{P}^{(2)}$ in $\mathrm{Pth}_{\boldsymbol{\mathcal{A}}^{(2)}, s}$ the $(2,[1])$-identity  second-order path on the $([1],2)$-target of $\mathfrak{P}^{(2)}$.

\textsf{(10)} For every $s\in S$, the partial binary operation of $1$-composition, i.e., $\circ_{s}^{1\mathbf{Pth}_{\boldsymbol{\mathcal{A}}^{(2)}}}$, is equal to the partial binary operation of $1$-composition of second-order paths introduced in Definition~\ref{DDPthComp}. We recall that, in Proposition~\ref{PDPthComp}, we proved that, given two second-order paths $\mathfrak{P}^{(2)}$ and $\mathfrak{Q}^{(2)}$ satisfying that 
$$
\mathrm{sc}^{([1],2)}_{s}
\left(
\mathfrak{Q}^{(2)}
\right)=
\mathrm{tg}^{([1],2)}_{s}\left(
\mathfrak{P}^{(2)}
\right),
$$
then the $1$-composition  $\mathfrak{Q}^{(2)}\circ_{s}^{1\mathbf{Pth}_{\boldsymbol{\mathcal{A}}^{(2)}}}\mathfrak{P}^{(2)}$ is well-defined.

This completes the definition of the partial $\Sigma^{\boldsymbol{\mathcal{A}}^{(2)}}$-algebra $\mathbf{Pth}_{\boldsymbol{\mathcal{A}}^{(2)}}$.
\end{proof}

\begin{remark} For the partial $\Sigma^{\boldsymbol{\mathcal{A}}^{(2)}}$-algebra $\mathbf{Pth}_{\boldsymbol{\mathcal{A}}^{(2)}}$, we denote by $\mathbf{Pth}^{(1,2)}_{\boldsymbol{\mathcal{A}}^{(2)}}$ the $\Sigma^{\boldsymbol{\mathcal{A}}}$-algebra $\mathbf{in}^{\Sigma,(1,2)}(\mathbf{Pth}_{\boldsymbol{\mathcal{A}}^{(2)}})$. We will call $\mathbf{Pth}^{(1,2)}_{\boldsymbol{\mathcal{A}}}$ the \emph{$\Sigma^{\boldsymbol{\mathcal{A}}}$-reduct} of the partial $\Sigma^{\boldsymbol{\mathcal{A}}^{(2)}}$-algebra $\mathbf{Pth}_{\boldsymbol{\mathcal{A}}^{(2)}}$. Note that this partial $\Sigma^{\boldsymbol{\mathcal{A}}}$-algebra coincides with the partial $\Sigma^{\boldsymbol{\mathcal{A}}}$-algebra introduced in Proposition~\ref{PDPthCatAlg}.

For the partial $\Sigma^{\boldsymbol{\mathcal{A}}^{(2)}}$-algebra $\mathbf{Pth}_{\boldsymbol{\mathcal{A}}^{(2)}}$, we denote by $\mathbf{Pth}^{(0,2)}_{\boldsymbol{\mathcal{A}}^{(2)}}$ the $\Sigma$-algebra $\mathbf{in}^{\Sigma,(0,2)}(\mathbf{Pth}_{\boldsymbol{\mathcal{A}}^{(2)}})$. We will call $\mathbf{Pth}^{(0,2)}_{\boldsymbol{\mathcal{A}}}$ the \emph{$\Sigma$-reduct} of the partial $\Sigma^{\boldsymbol{\mathcal{A}}^{(2)}}$-algebra $\mathbf{Pth}_{\boldsymbol{\mathcal{A}}^{(2)}}$. Note that this $\Sigma$-algebra coincides with the  $\Sigma$-algebra introduced in Proposition~\ref{PDPthCatAlg}.
\end{remark}

\section{
\texorpdfstring
{The second-order Curry-Howard mapping}
{The second-order Curry-Howard mapping}
}
The previous results, as will be seen below, will allow us to introduce a Curry-Howard approach to the definition of the many-sorted path term class transformation determined by the second-order specification $\boldsymbol{\mathcal{A}}^{(2)} = (\boldsymbol{\mathcal{A}},\mathcal{A}^{(2)})$. To do this, We define, by Artinian recursion, an $S$-sorted mapping from $\mathrm{Pth}_{\boldsymbol{\mathcal{A}}^{(2)}}$ to $\mathrm{T}_{\Sigma^{\boldsymbol{\mathcal{A}}^{(2)}}}(X)$, the underlying $S$-sorted set of $\mathbf{T}_{\Sigma^{\boldsymbol{\mathcal{A}}^{(2)}}}(X)$, the free algebra on the categorial signature determined by $\mathcal{A}^{(2)}$ on $X$. In this way, every second-order path in $\boldsymbol{\mathcal{A}}^{(2)}$ will be denoted by a term in $\mathrm{T}_{\Sigma^{\boldsymbol{\mathcal{A}}^{(2)}}}(X)$. To stress the new situation we will refer to this mapping as the second-order Curry-Howard mapping and we will denote it by $\mathrm{CH}^{(2)}$.

\begin{restatable}{definition}{DDCH}
\label{DDCH}  
\index{Curry-Howard!second-order!$\mathrm{CH}^{(2)}$}
The \emph{second-order Curry-Howard mapping} is the $S$-sorted mapping
$$
\textstyle
\mathrm{CH}^{(2)}\colon\mathrm{Pth}_{\boldsymbol{\mathcal{A}}^{(2)}}\mor\mathrm{T}_{\Sigma^{\boldsymbol{\mathcal{A}}^{(2)}}}(X)
$$
defined by Artinian recursion on $(\coprod\mathrm{Pth}_{\boldsymbol{\mathcal{A}}^{(2)}}, \leq_{\mathbf{Pth}_{\boldsymbol{\mathcal{A}}^{(2)}}})$ as follows.

\textsf{Base step of the Artinian recursion}.

Let $(\mathfrak{P}^{(2)},s)$ be a minimal element of $(\coprod\mathrm{Pth}_{\boldsymbol{\mathcal{A}}^{(2)}}, \leq_{\mathbf{Pth}_{\boldsymbol{\mathcal{A}}^{(2)}}})$. Then, by Proposition~\ref{PDMinimal}, the second-order path $\mathfrak{P}^{(2)}$ is either~(1) an $(2,[1])$-identity second-order path or~(2) a   second-order echelon.

If~(1), i.e., if $\mathfrak{P}^{(2)}$ is a $(2,[1])$-identity second-order path, then $\mathfrak{P}^{(2)}=\mathrm{ip}^{(2,[1])\sharp}_{s}([P]_{s})$ for some path term class $[P]_{s}$ in $[\mathrm{PT}_{\boldsymbol{\mathcal{A}}}]_{s}$. We define $\mathrm{CH}^{(2)}_{s}(\mathfrak{P}^{(2)})$ to be the term in $\mathrm{T}_{\Sigma^{\boldsymbol{\mathcal{A}}^{(2)}}}(X)_{s}$ given by the lift of the original Curry-Howard mapping applied to the interpretation as a path of its defining path term class, i.e.,
$$
\mathrm{CH}^{(2)}_{s}\left(
\mathfrak{P}^{(2)}
\right)=
\eta^{(2,1)\sharp}_{s}\left(
\mathrm{CH}^{(1)\mathrm{m}}_{s}\left(
\mathrm{ip}^{([1],X)@}_{s}\left(
\left[P
\right]_{s}
\right)\right)\right).
$$

If~(2), i.e., if $\mathfrak{P}^{(2)}$ is a   second-order echelon associated to a second-order rewrite rule $\mathfrak{p}^{(2)}=([M]_{s}, [N]_{s})$, that is, if $\mathfrak{P}^{(2)}$ has the form
\begin{center}
\begin{tikzpicture}
[ACliment/.style={-{To [angle'=45, length=5.75pt, width=4pt, round]},font=\scriptsize},
AClimentD/.style={double equal sign distance, -implies, font=\scriptsize}
]
\node[] (1) at (0,0) [] {$\mathfrak{P}^{(2)}\colon [M]_{s}$};
\node[] (2) at (5,0) [] {$[N]_{s}$.};
\draw[AClimentD]  (1) to node [above]
{$(\mathfrak{p}^{(2)},\mathrm{id}^{\mathrm{T}_{\Sigma^{\boldsymbol{\mathcal{A}}}}(X)_{s}})$} (2);
\end{tikzpicture}
\end{center}
then we define $\mathrm{CH}^{(2)}_{s}(\mathfrak{P}^{(2)})$ to be the syntactic representation of the unique second-order rewrite rule occurring in $\mathfrak{P}^{(2)}$, i.e., 
$$
\mathrm{CH}^{(2)}_{s}\left(
\mathfrak{P}^{(2)}
\right)=\mathfrak{p}^{(2)\mathbf{T}_{\Sigma^{\boldsymbol{\mathcal{A}}^{(2)}}}(X)}.
$$

\textsf{Inductive step of the Artinian recursion}.

Let $(\mathfrak{P}^{(2)},s)$ be a non-minimal element of $(\coprod\mathrm{Pth}_{\boldsymbol{\mathcal{A}}^{(2)}}, \leq_{\mathbf{Pth}_{\boldsymbol{\mathcal{A}}^{(2)}}})$.
We can assume that $\mathfrak{P}^{(2)}$ is not a $(2,[1])$-identity second-order path, since those second-order paths already have an image for the second-order Curry-Howard mapping. Let us suppose that, for every sort $t\in S$ and every second-order path $\mathfrak{Q}^{(2)}\in\mathrm{Pth}_{\boldsymbol{\mathcal{A}}^{(2)},t}$, if $(\mathfrak{Q}^{(2)},t)<_{\mathbf{Pth}_{\boldsymbol{\mathcal{A}}^{(2)}}}(\mathfrak{P}^{(2)},s)$, then the value of the second-order Curry-Howard mapping at $\mathfrak{Q}^{(2)}$, i.e., $\mathrm{CH}^{(2)}_{t}(\mathfrak{Q}^{(2)})$, has already been defined.

By Lemma~\ref{LDOrdI}, we have that $\mathfrak{P}^{(2)}$ is either~(1) a second-order path of length strictly greater than one containing at least one   second-order echelon or~(2) an echelonless second-order path.

If~(1), i.e., if $\mathfrak{P}^{(2)}$ is a second-order path of length strictly greater than one containing at least one   second-order echelon, then let $i\in \bb{\mathfrak{P}^{(2)}}$ be the first index for which the one-step subpath $\mathfrak{P}^{(2),i,i}$ of $\mathfrak{P}^{(2)}$ is a   second-order echelon. We consider different cases for $i$ according to the cases presented in Definition~\ref{DDOrd}.

If $i=0$, we have that the pairs $(\mathfrak{P}^{(2),0,0},s)$ and $(\mathfrak{P}^{(2),1,\bb{\mathfrak{P}^{(2)}}-1},s)$ $\prec_{\mathbf{Pth}_{\boldsymbol{\mathcal{A}}^{(2)}}}$-precede the pair $(\mathfrak{P}^{(2)},s)$. Therefore, the values of the second-order Curry-Howard mapping at $\mathfrak{P}^{(2),0,0}$ and $\mathfrak{P}^{(2),1,\bb{\mathfrak{P}^{(2)}}-1}$, respectively, have already been defined. 

In this case, we set
$$
\mathrm{CH}^{(2)}_{s}\left(\mathfrak{P}^{(2)}
\right)=
\mathrm{CH}^{(2)}_{s}
\left(
\mathfrak{P}^{(2),1,\bb{\mathfrak{P}^{(2)}}-1}
\right)
\circ_{s}^{1\mathbf{T}_{\Sigma^{\boldsymbol{\mathcal{A}}^{(2)}}}(X)}
\mathrm{CH}^{(2)}_{s}
\left(
\mathfrak{P}^{(2),0,0}
\right).
$$

If $i\neq 0$, we have that the pairs $(\mathfrak{P}^{(2),0,i-1},s)$ and $(\mathfrak{P}^{(2),i,\bb{\mathfrak{P}^{(2)}}-1},s)$ $\prec_{\mathbf{Pth}_{\boldsymbol{\mathcal{A}}^{(2)}}}$-precede the pair $(\mathfrak{P}^{(2)},s)$. Therefore, the values of the second-order Curry-Howard mapping at $\mathfrak{P}^{(2),0,i-1}$ and $\mathfrak{P}^{(2),i,\bb{\mathfrak{P}^{(2)}}-1}$, respectively, have already been defined. 

In this case, we set
$$
\mathrm{CH}^{(2)}_{s}\left(
\mathfrak{P}^{(2)}
\right)=
\mathrm{CH}^{(2)}_{s}
\left(\mathfrak{P}^{(2),i,\bb{\mathfrak{P}^{(2)}}-1}
\right)
\circ_{s}^{1\mathbf{T}_{\Sigma^{\boldsymbol{\mathcal{A}}^{(2)}}}(X)}
\mathrm{CH}^{(2)}_{s}\left(
\mathfrak{P}^{(2),0,i-1}
\right).
$$

This finishes the definition of the value of the second-order Curry-Howard mapping at a second-order path of length strictly greater than one containing at least one   second-order echelon.

If~(2), i.e., if $\mathfrak{P}^{(2)}$ is an echelonless second-order path in $\mathrm{Pth}_{\boldsymbol{\mathcal{A}}^{(2)},s}$. It could be the case that~(2.1) $\mathfrak{P}^{(2)}$ is not head-constant. Then let $i\in \bb{\mathfrak{P}^{(2)}}$ be the maximum index for which the subpath $\mathfrak{P}^{(2),0,i}$ of $\mathfrak{P}^{(2)}$ is a head-constant, echelonless second-order path. Note that the pairs $(\mathfrak{P}^{(2),0,i},s)$ and $(\mathfrak{P}^{(2),i+1,\bb{\mathfrak{P}^{(2)}}-1},s)$ $\prec_{\mathbf{Pth}_{\boldsymbol{\mathcal{A}}^{(2)}}}$-precede the pair $(\mathfrak{P}^{(2)},s)$. Therefore, the values of the second-order Curry-Howard mapping at $\mathfrak{P}^{(2),0,i}$ and $\mathfrak{P}^{(2),i+1,\bb{\mathfrak{P}^{(2)}}-1}$, respectively, have already been defined. 

In this case, we set
$$
\mathrm{CH}^{(2)}_{s}
\left(\mathfrak{P}^{(2)}
\right)=
\mathrm{CH}^{(2)}_{s}
\left(\mathfrak{P}^{(2),i+1,\bb{\mathfrak{P}^{(2)}}-1}
\right)
\circ_{s}^{1\mathbf{T}_{\Sigma^{\boldsymbol{\mathcal{A}}^{(2)}}}(X)}
\mathrm{CH}^{(2)}_{s}
\left(\mathfrak{P}^{(2),0,i}
\right).
$$

Therefore we are left with the case of $\mathfrak{P}^{(2)}$ being a head-constant echelonless second-order path. It could be the case that~(2.2)
 $\mathfrak{P}^{(2)}$ is not coherent. Then let $i\in \bb{\mathfrak{P}^{(2)}}$ be the maximum index for which the subpath $\mathfrak{P}^{(2),0,i}$ of $\mathfrak{P}^{(2)}$ is a coherent head-constant echelonless second-order path. Note that the pairs $(\mathfrak{P}^{(2),0,i},s)$ and $(\mathfrak{P}^{(2),i+1,\bb{\mathfrak{P}^{(2)}}-1},s)$ $\prec_{\mathbf{Pth}_{\boldsymbol{\mathcal{A}}^{(2)}}}$-precede the pair $(\mathfrak{P}^{(2)},s)$. Therefore, the values of the second-order Curry-Howard mapping at $\mathfrak{P}^{(2),0,i}$ and $\mathfrak{P}^{(2),i+1,\bb{\mathfrak{P}^{(2)}}-1}$, respectively, have already been defined. 

In this case, we set
$$
\mathrm{CH}^{(2)}_{s}
\left(
\mathfrak{P}^{(2)}
\right)=
\mathrm{CH}^{(2)}_{s}
\left(\mathfrak{P}^{(2),i+1,\bb{\mathfrak{P}^{(2)}}-1}
\right)
\circ_{s}^{1\mathbf{T}_{\Sigma^{\boldsymbol{\mathcal{A}}^{(2)}}}(X)}
\mathrm{CH}^{(2)}_{s}
\left(\mathfrak{P}^{(2),0,i}
\right).
$$

Therefore we are left with the case~(2.3) of $\mathfrak{P}^{(2)}$ being a coherent head-constant echelonless second-order path. Under this setting, the conditions for the second-order extraction algorithm, that is, Lemma~\ref{LDPthExtract} are fulfilled. Then there exists a unique word $\mathbf{s}\in S^{\star}-\{\lambda\}$ and a unique operation symbol $\tau\in \Sigma^{\boldsymbol{\mathcal{A}}}_{\mathbf{s},s}$ associated to $\mathfrak{P}^{(2)}$. Let $(\mathfrak{P}^{(2)}_{j})_{j\in\bb{\mathbf{s}}}$ be the family of second-order paths in $\mathrm{Pth}_{\boldsymbol{\mathcal{A}}^{(2)},\mathbf{s}}$ which, in virtue of Lemma~\ref{LDPthExtract}, we can extract from $\mathfrak{P}^{(2)}$. Note that, for every $j\in\bb{\mathbf{s}}$, we have that 
$(\mathfrak{P}^{(2)}_{j},s_{j})\prec_{\mathbf{Pth}_{\boldsymbol{\mathcal{A}}^{(2)}}}(\mathfrak{P}^{(2)},s)$.  Therefore, for every $j\in\bb{\mathbf{s}}$, the value of the second-order Curry-Howard mapping at $\mathfrak{P}^{(2)}_{j}$ has already been defined.

In this case, we set
$$
\mathrm{CH}^{(2)}_{s}\left(
\mathfrak{P}^{(2)}
\right)=
\tau^{\mathbf{T}_{\Sigma^{\boldsymbol{\mathcal{A}}^{(2)}}}(X)}
\left(
\left(
\mathrm{CH}^{(2)}_{s_{j}}
\left(\mathfrak{P}^{(2)}_{j}
\right)\right)_{j\in\bb{\mathbf{s}}}
\right).
$$

This finishes the definition of the value of the second-order Curry-Howard mapping at an echelonless second-order path.

This completes the definition of the second-order Curry-Howard mapping.
\end{restatable}

\begin{remark}\label{RDCH} The second-order Curry-Howard mapping associates to every second-order path of sort $s\in S$ a term in $\mathrm{T}_{\Sigma^{\boldsymbol{\mathcal{A}}^{(2)}}}(X)_{s}$.  In this regard, we highlight  how crucial has been the categorial second-order expansion of the original signature. The reader will realize that the second-order Curry-Howard mapping equates a syntactic transformation of a path term class --a second-order path on paths on terms-- with a term, at the expense of adding complexity on the term's side. This term contains the most important features of the second-order path and will later be used to reconstruct another second-order path, not necessarily equal to the original, but more standardized, where the rewrite rules acting in parallel are applied following different derivation strategies.
\end{remark}

\section{
\texorpdfstring
{The behaviour of the second-order Curry-Howard mapping}
{Behaviour}
}

The following propositions will provide a deeper understanding of the just defined second-order Curry-Howard mapping. We first study how the second-order Curry-Howard mapping is correlated with the mappings arising from $X$.

\begin{figure}
\begin{center}
\begin{tikzpicture}
[ACliment/.style={-{To [angle'=45, length=5.75pt, width=4pt, round]}},scale=0.8]
\node[] (x) at (-4,-2) [] {$X$};
\node[] (pth12) at (2,-2) [] {$\mathrm{Pth}_{\boldsymbol{\mathcal{A}}}$};
\node[] (pt12) at (2,-4) [] {$\mathrm{PT}_{\boldsymbol{\mathcal{A}}}$};
\node[] (pth1) at (8,-2) [] {$[\mathrm{Pth}_{\boldsymbol{\mathcal{A}}}]$};
\node[] (pt1) at (8,-4) [] {$[\mathrm{PT}_{\boldsymbol{\mathcal{A}}}]$};
\node[] (pth2) at (8,-6) [] {$\mathrm{Pth}_{\boldsymbol{\mathcal{A}}^{(2)}}$};
\node[] (t2) at (8,-8) [] {$\mathrm{T}_{\Sigma^{\boldsymbol{\mathcal{A}}^{(2)}}}(X)$};

\draw[ACliment]  (x) to node [above right]
{$\textstyle {\mathrm{ip}}^{(1,X)}$} (pth12);
\draw[ACliment, bend right=10]  (x) to node  [above right]
{$\textstyle {\eta}^{(1,X)}$} (pt12);
\draw[ACliment, bend right=20]  (x) to node  [midway, fill=white, pos=0.8]
{$\textstyle \mathrm{ip}^{(2,X)}$} (pth2);
\draw[ACliment, bend right=30]  (pt12) to node  [midway, fill=white]
{$\textstyle \eta^{(2,1)\sharp}$} (t2);
\draw[ACliment, bend right]  (x) to node [ below left, pos=0.8]
{$\textstyle \eta^{(2,X)}$} (t2);

\draw[ACliment]  (pth12) to node [above]
{$\mathrm{pr}^{\mathrm{Ker}(\mathrm{CH}^{(1)})}$} (pth1);

\draw[ACliment]  (pt12) to node [below]
{$\mathrm{pr}^{\Theta^{[1]}}$} (pt1);

\draw[ACliment]  (pth1) to node [midway, fill=white]
{$\mathrm{CH}^{(1)\mathrm{m}}$} (pt12);

\draw[ACliment] 
($(pt1)+(0,-.35)$) to node [right] {$\textstyle \mathrm{ip}^{(2,[1])\sharp}$}  ($(pth2)+(0,.35)$);

\draw[ACliment] 
($(pth1)+(-.15,-.35)$) to node [left] {$\textstyle \mathrm{CH}^{[1]}$}  ($(pt1)+(-.15,.35)$);
\draw[ACliment] 
($(pt1)+(.15,+.35)$) to node [right] {$\textstyle \mathrm{ip}^{([1],X)@}$}  ($(pth1)+(.15,-.35)$);

\draw[ACliment] 
($(pth12)+(-.15,-.35)$) to node [left] {$\textstyle \mathrm{CH}^{(1)}$}  ($(pt12)+(-.15,.35)$);
\draw[ACliment] 
($(pt12)+(.15,+.35)$) to node [right] {$\textstyle \mathrm{ip}^{(1,X)@}$}  ($(pth12)+(.15,-.35)$);

\draw[ACliment] 
($(pth2)+(-.15,-.35)$) to node [left] {$\textstyle \mathrm{CH}^{(2)}$}  ($(t2)+(-.15,.35)$);
\end{tikzpicture}
\end{center}
\caption{Behaviour of $\mathrm{CH}^{(2)}$ relative to $X$ at layers 1 \& 2.}
\label{FDCHDUId}
\end{figure}

Next proposition describes the image under the second-order Curry-Howard mapping of $(2,[1])$-identity second-order paths.

\begin{proposition}\label{PDCHDUId}
For the second-order Curry-Howard mapping we have that
\allowdisplaybreaks
\begin{align*}
\mathrm{(i)}\, \mathrm{CH}^{(2)}\circ\mathrm{ip}^{(2,X)}=\eta^{(2,X)}
&&
\mathrm{(ii)}\, \mathrm{CH}^{(2)}\circ\mathrm{ip}^{(2,[1])\sharp}=\eta^{(2,1)\sharp}\circ\mathrm{CH}^{(1)\mathrm{m}}\circ\mathrm{ip}^{([1],X)@}.
\end{align*}
\end{proposition}
\begin{proof}
The reader is advised to consult the diagram presented in Figure~\ref{FDCHDUId}.

Let $s$ be a sort in $S$ and let $x$ be a variable in $X_{s}$, then the following chain of equalities holds 
\allowdisplaybreaks
\begin{align*}
\mathrm{CH}^{(2)}_{s}
\left(
\mathrm{ip}^{(2,X)}_{s}
\left(x
\right)\right)&=
\mathrm{CH}^{(2)}_{s}
\left(
\left[
x
\right]_{s},\lambda,\lambda
\right)
\tag{1}
\\&=
\eta^{(2,1)\sharp}_{s}\left(
\mathrm{CH}^{(1)\mathrm{m}}_{s}\left(
\mathrm{ip}^{([1],X)@}_{s}\left(
\left[x
\right]_{s}
\right)\right)\right)
\tag{2}
\\&=
\eta^{(2,1)\sharp}_{s}\left(
\mathrm{CH}^{(1)\mathrm{m}}_{s}\left(
\mathrm{ip}^{([1],X)@}_{s}\left(
\eta^{([1],X)}_{s}\left(
x
\right)\right)\right)\right)
\tag{3}
\\&=
\eta^{(2,1)\sharp}_{s}\left(
\mathrm{CH}^{(1)\mathrm{m}}_{s}\left(
\mathrm{ip}^{([1],X)}_{s}\left(
x
\right)\right)\right)
\tag{4}
\\&=
\eta^{(2,1)\sharp}_{s}\left(
\mathrm{CH}^{(1)}_{s}\left(
\mathrm{ip}^{(1,X)}_{s}\left(
x
\right)\right)\right)
\tag{5}
\\&=
\eta^{(2,1)\sharp}_{s}\left(
\eta^{(1,X)}_{s}\left(
x
\right)\right)
\tag{6}
\\&=
\eta^{(2,X)}_{s}\left(x
\right).\tag{6}
\end{align*}

The first equality unravels the definition of the mapping $\mathrm{ip}^{(2,X)}$; the second equality applies the second-order Curry-Howard mapping at a $(2,[1])$-identity second-order path; the third equality follows from Proposition~\ref{PPTQXEq}, since $\mathrm{ip}^{([1],X)@}\circ\eta^{([1],X)}=\mathrm{ip}^{([1],X)}$; the fourth equality unravels the definition of the mappings $\mathrm{CH}^{(1)\mathrm{m}}$ and $\mathrm{ip}^{([1],X)@}$; the fifth equality  follows from Proposition~\ref{PCHId}, since $\mathrm{CH}^{(1)}\circ\mathrm{ip}^{(1,X)}=\mathrm{\eta}^{(1,X)}$; finally, the last equality follows from Proposition~\ref{PDEmb}, since $\eta^{(2,1)\sharp}\circ\eta^{(1,X)}=\eta^{(2,X)}$.

Now, let $s$ be a sort in $S$ and let $[P]_{s}$ be a path term class in $[\mathrm{PT}_{\boldsymbol{\mathcal{A}}}]_{s}$, then, according to Definition~\ref{DDCH},  we have that
\allowdisplaybreaks
\begin{align*}
\mathrm{CH}^{(2)}_{s}
\left(\mathrm{ip}^{(2,[1])\sharp}_{s}
\left(\left[
P
\right]_{s}
\right)\right)
&=
\eta^{(2,1)\sharp}_{s}
\left(\mathrm{CH}^{(1)\mathrm{m}}_{s}\left(
\mathrm{ip}^{([1],X)@}_{s}
\left(\left[
P
\right]_{s}
\right)\right)\right).
\end{align*}

This finishes the proof.
\end{proof}

We next prove that the second-order Curry-Howard mapping is a $\Sigma$-homomorphism from the $\Sigma$-algebra $\mathbf{Pth}^{(0,2)}_{\boldsymbol{\mathcal{A}}^{(2)}}$ to the $\Sigma$-algebra $\mathbf{T}^{(0,2)}_{\Sigma^{\boldsymbol{\mathcal{A}}^{(2)}}}(X)$.

\begin{restatable}{proposition}{PDCHHom}
\label{PDCHHom}
The second-order Curry-Howard mapping is a $\Sigma$-homomorphism from $\mathbf{Pth}^{(0,2)}_{\boldsymbol{\mathcal{A}}^{(2)}}$ to $\mathbf{T}^{(0,2)}_{\Sigma^{\boldsymbol{\mathcal{A}}^{(2)}}}(X)$. 
\end{restatable}

\begin{proof}
Let $(\mathbf{s},s)$ be an element of  $S^{\star}\times S$ and $\sigma$ an operation symbol in $\Sigma_{\mathbf{s},s}$. 

If $\mathbf{s}=\lambda$, then the following chain of equalities holds
\begin{align*}
\mathrm{CH}^{(2)}_{s}
\left(
\sigma^{\mathbf{Pth}_{\boldsymbol{\mathcal{A}}^{(2)}}}
\right)
&=
\mathrm{CH}^{(2)}_{s}
\left(
\mathrm{ip}^{(2,[1])\sharp}_{s}
\left(
\left[
\sigma^{\mathbf{PT}_{\boldsymbol{\mathcal{A}}}}
\right]_{s}
\right)\right)
\tag{1}
\\&=
\eta^{(2,1)\sharp}_{s}
\left(
\mathrm{CH}^{(1)\mathrm{m}}_{s}
\left(
\mathrm{ip}^{([1],X)@}_{s}
\left(
\left[
\sigma^{\mathbf{PT}_{\boldsymbol{\mathcal{A}}}}
\right]_{s}
\right)\right)\right)
\tag{2}
\\&=
\eta^{(2,1)\sharp}_{s}
\left(
\mathrm{CH}^{(1)\mathrm{m}}_{s}
\left(
\left[
\sigma^{\mathbf{Pth}_{\boldsymbol{\mathcal{A}}}}
\right]_{s}
\right)\right)
\tag{3}
\\&=
\eta^{(2,1)\sharp}_{s}
\left(
\sigma^{\mathbf{PT}_{\boldsymbol{\mathcal{A}}}}
\right)
\tag{4}
\\&=
\sigma^{\mathbf{T}_{\Sigma^{\boldsymbol{\mathcal{A}}^{(2)}}}(X)}.
\tag{5}
\end{align*}

The first equality follows from the fact that, for a constant $\sigma$ in $\Sigma_{\lambda,s}$, then, following Remark~\ref{RDConsSigma}, the interpretation of $\sigma$ in $\mathbf{Pth}_{\boldsymbol{\mathcal{A}}^{(2)}}$ is given by $
\sigma^{\mathbf{Pth}_{\boldsymbol{\mathcal{A}}^{(2)}}}=\mathrm{ip}^{(2,[1])\sharp}_{s}(
[\sigma^{\mathbf{PT}_{\boldsymbol{\mathcal{A}}}}]_{s}
)$; the second equality follows from Proposition~\ref{PDCHDUId}; the third equality follows from the fact that, by Theorem~\ref{TIso}, $\mathrm{ip}^{([1],X)@}$ is a $\Sigma^{\boldsymbol{\mathcal{A}}}$-homomorphism; the fourth equality applies the monomorphic Curry-Howard mapping at a path class; finally, the last equality follows from the fact that $\eta^{(2,1)\sharp}$ is a $\Sigma^{\boldsymbol{\mathcal{A}}}$-homomorphism from Proposition~\ref{PDUEmb}.

We now consider the case in which $\mathbf{s}\neq \lambda$. Let $(\mathfrak{P}^{(2)}_{j})_{j\in\bb{\mathbf{s}}}$ be a family of second-order paths in $\mathrm{Pth}_{\boldsymbol{\mathcal{A}}^{(2)},\mathbf{s}}$. We consider different cases according to the nature of the family $(\mathfrak{P}^{(2)}_{j})_{j\in\bb{\mathbf{s}}}$. It could be the case that either (1), for every $j\in \bb{\mathbf{s}}$, $\mathfrak{P}^{(2)}_{j}$ is a $(2,[1])$-identity second-order path or (2), there exists an index $j\in\bb{\mathbf{s}}$ for which $\mathfrak{P}^{(2)}_{j}$ is a non-$(2,[1])$-identity second-order path.

If~(1), then for every $j\in\bb{\mathbf{s}}$, $\mathfrak{P}^{(2)}_{j}$ is equal to $\mathrm{ip}^{(2,[1])\sharp}_{s_{j}}([P_{j}]_{s_{j}})$, for some path term class $[P_{j}]_{s_{j}}$ in $[\mathrm{PT}_{\boldsymbol{\mathcal{A}}}]_{s_{j}}$.  In this case, the following chain of equalities holds
\begin{flushleft}
$
\mathrm{CH}^{(2)}_{s}\left(
\sigma^{\mathbf{Pth}_{\boldsymbol{\mathcal{A}}^{(2)}}}
\left(\left(
\mathfrak{P}^{(2)}_{j}
\right)_{j\in\bb{\mathbf{s}}}
\right)\right)$
\allowdisplaybreaks
\begin{align*}
\,
&=
\mathrm{CH}^{(2)}_{s}\left(
\sigma^{\mathbf{Pth}_{\boldsymbol{\mathcal{A}}^{(2)}}}
\left(\left(
\mathrm{ip}^{(2,[1])\sharp}_{s_{j}}\left(
\left[
P_{j}
\right]_{s_{j}}
\right)\right)_{j\in\bb{\mathbf{s}}}
\right)\right)
\tag{1}
\\&=
\mathrm{CH}^{(2)}_{s}\left(
\mathrm{ip}^{(2,[1])\sharp}_{s}\left(
\sigma
^{[\mathbf{PT}_{\boldsymbol{\mathcal{A}}}]}
\left(\left(
\left[P_{j}
\right]_{s_{j}}
\right)_{j\in\bb{\mathbf{s}}}
\right)\right)\right)
\tag{2}
\\&=
\eta^{(2,1)\sharp}_{s}\left(
\mathrm{CH}^{(1)\mathrm{m}}_{s}\left(
\mathrm{ip}^{([1],X)@}_{s}\left(
\sigma
^{[\mathbf{PT}_{\boldsymbol{\mathcal{A}}}]}
\left(\left(
\left[P_{j}
\right]_{s_{j}}
\right)_{j\in\bb{\mathbf{s}}}
\right)\right)\right)\right)
\tag{3}
\\&=
\eta^{(2,1)\sharp}_{s}\left(
\mathrm{CH}^{(1)\mathrm{m}}_{s}\left(
\sigma^{[\mathbf{Pth}_{\boldsymbol{\mathcal{A}}}]}
\left(\left(
\mathrm{ip}^{([1],X)@}_{s_{j}}
\left(
\left[
P_{j}
\right]_{s_{j}}
\right)\right)_{j\in\bb{\mathbf{s}}}
\right)\right)\right)
\tag{4}
\\&=
\eta^{(2,1)\sharp}_{s}\left(
\mathrm{CH}^{(1)\mathrm{m}}_{s}\left(
\sigma^{[\mathbf{Pth}_{\boldsymbol{\mathcal{A}}}]}
\left(\left(
\left[\mathrm{ip}^{([1],X)@}_{s_{j}}
\left(P_{j}\right)
\right]_{s_{j}}
\right)_{j\in\bb{\mathbf{s}}}
\right)\right)\right)
\tag{5}
\\&=
\eta^{(2,1)\sharp}_{s}\left(
\mathrm{CH}^{(1)\mathrm{m}}_{s}\left(
\left[
\sigma^{\mathbf{Pth}_{\boldsymbol{\mathcal{A}}}}
\left(\left(
\mathrm{ip}^{(1,X)@}_{s_{j}}
\left(
P_{j}
\right)\right)_{j\in\bb{\mathbf{s}}}
\right)\right]_{s}
\right)\right)
\tag{6}
\\&=
\eta^{(2,1)\sharp}_{s}\left(
\mathrm{CH}^{(1)}_{s}\left(
\sigma^{\mathbf{Pth}_{\boldsymbol{\mathcal{A}}}}
\left(\left(
\mathrm{ip}^{(1,X)@}_{s_{j}}\left(
P_{j}
\right)\right)_{j\in\bb{\mathbf{s}}}
\right)\right)\right)
\tag{7}
\\&=
\eta^{(2,1)\sharp}_{s}\left(
\sigma^{\mathbf{PT}_{\boldsymbol{\mathcal{A}}}}
\left(\left(
\mathrm{CH}^{(1)}_{s_{j}}
\left(
\mathrm{ip}^{(1,X)@}_{s_{j}}\left(
P_{j}
\right)\right)\right)_{j\in\bb{\mathbf{s}}}
\right)\right)
\tag{8}
\\&=
\sigma^{\mathbf{T}_{\Sigma^{\boldsymbol{\mathcal{A}}^{(2)}}}(X)}
\left(\left(
\eta^{(2,1)\sharp}_{s_{j}}\left(
\mathrm{CH}^{(1)}_{s_{j}}\left(
\mathrm{ip}^{(1,X)@}_{s_{j}}\left(
P_{j}
\right)\right)\right)\right)_{j\in\bb{\mathbf{s}}}
\right)
\tag{9}
\\&=
\sigma^{\mathbf{T}_{\Sigma^{\boldsymbol{\mathcal{A}}^{(2)}}}(X)}
\left(\left(
\eta^{(2,1)\sharp}_{s_{j}}\left(
\mathrm{CH}^{(1)\mathrm{m}}_{s_{j}}\left(
\left[\mathrm{ip}^{(1,X)@}_{s_{j}}
\left(
P_{j}
\right)\right]_{s_{j}}
\right)\right)\right)_{j\in\bb{\mathbf{s}}}
\right)
\tag{10}
\\&=
\sigma^{\mathbf{T}_{\Sigma^{\boldsymbol{\mathcal{A}}^{(2)}}}(X)}
\left(\left(
\eta^{(2,1)\sharp}_{s_{j}}\left(
\mathrm{CH}^{(1)\mathrm{m}}_{s_{j}}\left(
\mathrm{ip}^{([1],X)@}_{s_{j}}\left(
\left[P_{j}
\right]_{s_{j}}
\right)\right)\right)\right)_{j\in\bb{\mathbf{s}}}
\right)
\tag{11}
\\&=
\sigma^{\mathbf{T}_{\Sigma^{\boldsymbol{\mathcal{A}}^{(2)}}}(X)}
\left(\left(
\mathrm{CH}^{(2)}_{s_{j}}\left(
\mathrm{ip}^{(2,[1])\sharp}_{s_{j}}\left(
\left[P_{j}
\right]_{s_{j}}
\right)\right)\right)_{j\in\bb{\mathbf{s}}}
\right)
\tag{12}
\\&=
\sigma^{\mathbf{T}_{\Sigma^{\boldsymbol{\mathcal{A}}^{(2)}}}(X)}
\left(\left(
\mathrm{CH}^{(2)}_{s_{j}}\left(
\mathfrak{P}^{(2)}_{j}
\right)\right)_{j\in\bb{\mathbf{s}}}
\right).
\tag{13}
\end{align*}
\end{flushleft}

The first equality follows from the fact that, for every $j\in\bb{\mathbf{s}}$, $\mathfrak{P}^{(2)}_{j}$ is a $(2,[1])$-identity second-order path; the second equality follows from the fact that, by Proposition~\ref{PDUIpCatHom}, $\mathrm{ip}^{(2,[1])\sharp}$ is a $\Sigma^{\boldsymbol{\mathcal{A}}}$-homomorphism from $[\mathbf{PT}_{\boldsymbol{\mathcal{A}}}]$ to $\mathbf{Pth}^{(1,2)}_{\boldsymbol{\mathcal{A}}^{(2)}}$; the third equality follows  follows from the characterization of the value for the second-order Curry-Howard mapping of $(2,[1])$-identity second-order paths introduced in Proposition~\ref{PDCHDUId};
the fourth equality follows from the fact that, by Theorem~\ref{TIso}, $\mathrm{ip}^{([1],X)@}$ is a $\Sigma^{\boldsymbol{\mathcal{A}}}$-homomorphism from $[\mathbf{PT}_{\boldsymbol{\mathcal{A}}}]$ to $[\mathbf{Pth}_{\boldsymbol{\mathcal{A}}}]$; the fifth equality applies the mapping $\mathrm{ip}^{([1],X)@}$ at each path term class; the sixth equality recovers the interpretation of the operation symbol $\sigma$ in the many-sorted partial $\Sigma^{\boldsymbol{\mathcal{A}}}$ of path classes $[\mathbf{Pth}_{\boldsymbol{\mathcal{A}}}]$; the seventh equality applies the mapping $\mathrm{CH}^{(1)\mathrm{m}}$ at the given path class; the eighth equality follows from the fact that, by Proposition~\ref{PCHHom}, $\mathrm{CH}^{(1)}$ is a $\Sigma$-homomorphism from $\mathbf{Pth}_{\boldsymbol{\mathcal{A}}}^{(0,1)}$ to $\mathbf{T}_{\Sigma^{\boldsymbol{\mathcal{A}}}}^{(0,1)}(X)$ and from the fact that, by Corollary~\ref{CPTCH}, the Curry-Howard mapping corestricts to $\mathrm{PT}_{\boldsymbol{\mathcal{A}}}$; the ninth equality follows from the fact that, by Proposition~\ref{PDUEmb}, $\eta^{(2,1)\sharp}$ is a $\Sigma^{\boldsymbol{\mathcal{A}}}$ from $\mathbf{PT}_{\boldsymbol{\mathcal{A}}}$ to $\mathbf{T}_{\Sigma^{\boldsymbol{\mathcal{A}}^{(2)}}}^{(1,2)}(X)$; the tenth equality recovers the definition of the mapping $\mathrm{CH}^{(1)\mathrm{m}}$; the eleventh equality recovers the definition of the mapping $\mathrm{ip}^{([1],X)@}$; the twelfth equality follows from the characterization of the value for the second-order Curry-Howard mapping of $(2,[1])$-identity second-order paths introduced in Proposition~\ref{PDCHDUId}; finally the last equality recovers, for every $j\in\bb{\mathbf{s}}$, the definition of  $\mathfrak{P}^{(2)}_{j}$ as a $(2,[1])$-identity second-order path.

This proves Case~(1).

If~(2), i.e., if there exists some index $j\in\bb{\mathbf{s}}$ for which $\mathfrak{P}^{(2)}_{j}$ is a non-$(2,[1])$-identity second-order path then, according to Corollary~\ref{CDPthWB}, $\sigma^{\mathbf{Pth}_{\boldsymbol{\mathcal{A}}^{(2)}}}((\mathfrak{P}^{(2)}_{j})_{j\in\bb{\mathbf{s}}})$ is an echelonless, head-constant  and coherent second-order path. Moreover, according to Proposition~\ref{PDRecov}, the second-order path extraction algorithm from Lemma~\ref{LDPthExtract} applied to it retrieves the original family $(\mathfrak{P}^{(2)}_{j})_{j\in\bb{\mathbf{s}}}$.Then, following Definition~\ref{DDCH},  the second-order Curry-Howard mapping applied to $\sigma^{\mathbf{Pth}_{\boldsymbol{\mathcal{A}}^{(2)}}}((\mathfrak{P}^{(2)}_{j})_{j\in\bb{\mathbf{s}}})$ is given by
\begin{align*}
\mathrm{CH}^{(2)}_{s}\left(
\sigma^{\mathbf{Pth}_{\boldsymbol{\mathcal{A}}^{(2)}}}
\left(\left(\mathfrak{P}^{(2)}_{j}
\right)_{j\in\bb{\mathbf{s}}}
\right)\right)
=
\sigma^{\mathbf{T}_{\Sigma^{\boldsymbol{\mathcal{A}}^{(2)}}}(X)}
\left(\left(\mathrm{CH}^{(2)}_{s_{j}}\left(\mathfrak{P}^{(2)}_{j}
\right)\right)_{j\in\bb{\mathbf{s}}}
\right).
\end{align*}

Case~(2) follows.

This finishes the proof.
\end{proof}

Next proposition describes the image under the second-order Curry-Howard mapping of $(2,0)$-identity second-order paths.

\begin{figure}
\begin{center}
\begin{scaletikzpicturetowidth}{\textwidth}
\begin{tikzpicture}
[ACliment/.style={-{To [angle'=45, length=5.75pt, width=4pt, round]}},scale=0.8]
\node[] (x) at (-4,0) [] {$X$};
\node[] (t22) at (2,0) [] {$\mathrm{T}_{\Sigma}(X)$};
\node[] (pth12) at (2,-2) [] {$\mathrm{Pth}_{\boldsymbol{\mathcal{A}}}$};
\node[] (pt12) at (2,-4) [] {$\mathrm{PT}_{\boldsymbol{\mathcal{A}}}$};
\node[] (t) at (8,0) [] {$\mathrm{T}_{\Sigma}(X)$};
\node[] (pth1) at (8,-2) [] {$[\mathrm{Pth}_{\boldsymbol{\mathcal{A}}}]$};
\node[] (pt1) at (8,-4) [] {$[\mathrm{PT}_{\boldsymbol{\mathcal{A}}}]$};
\node[] (pth2) at (8,-6) [] {$\mathrm{Pth}_{\boldsymbol{\mathcal{A}}^{(2)}}$};
\node[] (t2) at (8,-8) [] {$\mathrm{T}_{\Sigma^{\boldsymbol{\mathcal{A}}^{(2)}}}(X)$};

\draw[ACliment]  (x) to node [above right]
{$\textstyle \eta^{(0,X)}$} (t22);
\draw[ACliment, bend right=10]   (x) to node [above right]
{$\textstyle {\mathrm{ip}}^{(1,X)}$} (pth12);
\draw[ACliment, bend right=20]  (x.320) to node  [above right]
{$\textstyle {\eta}^{(1,X)}$} (pt12);
\draw[ACliment, bend right=30]  (x) to node  [midway, fill=white, pos=0.8]
{$\textstyle \mathrm{ip}^{(2,X)}$} (pth2);

\draw[ACliment, bend right=30]  (pt12) to node  [midway, fill=white]
{$\textstyle \eta^{(2,1)\sharp}$} (t2.170);

\draw[ACliment, bend right]  (x) to node [ below left, pos=0.8]
{$\textstyle \eta^{(2,X)}$} (t2);

\draw[double equal sign distance]  (t) to node []
{} (t22);

\draw[ACliment]  (t22) to node [right]
{$\mathrm{ip}^{(1,0)\sharp}$} (pth12);

\draw[ACliment]  (pth12) to node [above]
{$\mathrm{pr}^{\mathrm{Ker}(\mathrm{CH}^{(1)})}$} (pth1);

\draw[ACliment]  (pt12) to node [below]
{$\mathrm{pr}^{\Theta^{[1]}}$} (pt1);

\draw[ACliment]  (pth1) to node [midway, fill=white]
{$\mathrm{CH}^{(1)\mathrm{m}}$} (pt12);

\draw[ACliment] 
($(t)+(0,-.35)$) to node [right] {$\textstyle \mathrm{ip}^{([1],0)\sharp}$}  ($(pth1)+(0,.35)$);
\draw[ACliment] 
($(pt1)+(0,-.35)$) to node [right] {$\textstyle \mathrm{ip}^{(2,[1])\sharp}$}  ($(pth2)+(0,.35)$);

\draw[ACliment] 
($(pth1)+(-.15,-.35)$) to node [left] {$\textstyle \mathrm{CH}^{[1]}$}  ($(pt1)+(-.15,.35)$);
\draw[ACliment] 
($(pt1)+(.15,+.35)$) to node [right] {$\textstyle \mathrm{ip}^{([1],X)@}$}  ($(pth1)+(.15,-.35)$);

\draw[ACliment] 
($(pth12)+(-.15,-.35)$) to node [left] {$\textstyle \mathrm{CH}^{(1)}$}  ($(pt12)+(-.15,.35)$);
\draw[ACliment] 
($(pt12)+(.15,+.35)$) to node [right] {$\textstyle \mathrm{ip}^{(1,X)@}$}  ($(pth12)+(.15,-.35)$);

\draw[ACliment] 
($(pth2)+(-.15,-.35)$) to node [left] {$\textstyle \mathrm{CH}^{(2)}$}  ($(t2)+(-.15,.35)$);

\draw[ACliment, rounded corners]
($(t.east)+(0,-.05)$) -- 
($(t)+(2.15,-.05)$) -- node [right] {$\textstyle \mathrm{ip}^{(2,0)\sharp}$} 
($(pth2)+(2.15,0)$) -- (pth2.east)
;

\draw[ACliment, rounded corners]
($(t.east)+(0,+.05)$) -- 
($(t)+(3.75,+.05)$) -- 
($(t2)+(3.75,0)$) -- node [above] {$\textstyle \eta^{(2,0)\sharp}$} 
 (t2.east)
;
\end{tikzpicture}
\end{scaletikzpicturetowidth}
\end{center}
\caption{Behaviour of $\mathrm{CH}^{(2)}$ relative to $X$ at layers 0 \& 2.}
\label{FDCHDZId}
\end{figure}

\begin{proposition}\label{PDCHDZId}
For the second-order Curry-Howard mapping we have that 
\begin{itemize}
\item[(i)] $
\mathrm{CH}^{(2)}\circ\mathrm{ip}^{(2,0)\sharp}=\eta^{(2,0)\sharp}.$
\end{itemize}
\end{proposition}
\begin{proof}
The reader is advised to consult the diagram presented in Figure~\ref{FDCHDZId}.

Let us recall that the composition $\mathrm{CH}^{(2)}\circ\mathrm{ip}^{(2,0)\sharp}$ is a $\Sigma$-homomorphism from $\mathbf{T}_{\Sigma}(X)$ to $\mathbf{T}_{\Sigma^{\boldsymbol{\mathcal{A}}^{(2)}}}^{(0,2)}(X)$ in virtue of Propositions~\ref{PDZHomIp} and~\ref{PDCHHom}.

Consequently, in virtue of the Universal Property of the many-sorted free $\Sigma$-algebra $\mathbf{T}_{\Sigma}(X)$, it suffices to prove that 
$$
\mathrm{CH}^{(2)}\circ\mathrm{ip}^{(2,0)\sharp}
\circ \eta^{(0,X)}=\eta^{(2,X)}
$$
to obtain the desired result.

Let $s$ be a sort in $S$ and let $x$ be a variable in $X_{s}$. Let us note that the following chain of equalities holds
\begin{flushleft}
$
\mathrm{CH}^{(2)}_{s}\left(
\mathrm{ip}^{(2,0)\sharp}_{s}\left(
\eta^{(0,X)}_{s}\left(
x
\right)\right)\right)
$
\allowdisplaybreaks
\begin{align*}
\qquad
&=
\mathrm{CH}^{(2)}_{s}\left(
\mathrm{ip}^{(2,[1])\sharp}_{s}\left(
\mathrm{CH}^{[1]}_{s}\left(
\mathrm{ip}^{([1],0)\sharp}_{s}\left(
\eta^{(0,X)}_{s}\left(
x
\right)\right)\right)\right)\right)
\tag{1}
\\&=
\mathrm{CH}^{(2)}_{s}\left(
\mathrm{ip}^{(2,[1])\sharp}_{s}\left(
\mathrm{CH}^{[1]}_{s}\left(
\mathrm{ip}^{([1],X)}_{s}\left(
x
\right)\right)\right)\right)
\tag{2}
\\&=
\mathrm{CH}^{(2)}_{s}\left(
\mathrm{ip}^{(2,[1])\sharp}_{s}\left(
\mathrm{CH}^{[1]}_{s}\left(
\left[\mathrm{ip}^{(1,X)}_{s}(x)
\right]_{s}
\right)\right)\right)
\tag{3}
\\&=
\mathrm{CH}^{(2)}_{s}\left(
\mathrm{ip}^{(2,[1])\sharp}_{s}\left(
\left[
\mathrm{CH}^{(1)}_{s}\left(
\mathrm{ip}^{(1,X)}_{s}\left(
x
\right)\right)
\right]_{s}
\right)\right)
\tag{4}
\\&=
\mathrm{CH}^{(2)}_{s}\left(
\mathrm{ip}^{(2,[1])\sharp}_{s}\left(
\left[
\eta^{(1,X)}_{s}\left(
x
\right)\right]_{s}
\right)\right)
\tag{5}
\\&=
\eta^{(2,1)\sharp}_{s}\left(
\mathrm{CH}^{(1)\mathrm{m}}_{s}\left(
\mathrm{ip}^{([1],X)@}_{s}\left(
\left[\eta^{(1,X)}_{s}
\left(
x
\right)\right]_{s}
\right)\right)\right)
\tag{6}
\\&=
\eta^{(2,1)\sharp}_{s}\left(
\mathrm{CH}^{(1)\mathrm{m}}_{s}\left(
\left[\mathrm{ip}^{(1,X)@}_{s}\left(
\eta^{(1,X)}_{s}(x)\right)\right]_{s}
\right)\right)
\tag{7}
\\&=
\eta^{(2,1)\sharp}_{s}\left(
\mathrm{CH}^{(1)\mathrm{m}}_{s}\left(
\left[\mathrm{ip}^{(1,X)}_{s}(x)
\right]_{s}
\right)\right)
\tag{8}
\\&=
\eta^{(2,1)\sharp}_{s}\left(
\mathrm{CH}^{(1)}_{s}\left(
\mathrm{ip}^{(1,X)}_{s}\left(
x
\right)\right)\right)
\tag{9}
\\&=
\eta^{(2,1)\sharp}_{s}\left(
\eta^{(1,X)}_{s}\left(
x
\right)\right)
\tag{10}
\\&=
\eta^{(2,X)}_{s}(x).
\tag{11}
\end{align*}
\end{flushleft} 

The first equality recovers the definition of the mapping $\mathrm{ip}^{(2,0)\sharp}$, presented in Definition~\ref{DDScTgZ}; the second equality follows from the fact that, by Proposition~\ref{PCHBasicEq}, we have that $\mathrm{ip}^{([1],0)\sharp}\circ\eta^{(0,X)}=\mathrm{ip}^{([1],X)}$; the third equality applies the mapping $\mathrm{ip}^{([1],X)}$, introduced in Definition~\ref{DCHEch}; the fourth equality applies the mapping $\mathrm{CH}^{[1]}$, introduced in Definition~\ref{DPTQCH}; the fifth equality follows from Proposition~\ref{PCHId}, since $\mathrm{CH}^{(1)}\circ\mathrm{ip}^{(1,X)}=\eta^{(1,X)}$; the sixth equality follows from the characterization of the value for the second-order Curry-Howard mapping of $(2,[1])$-identity second-order paths introduced in  Proposition~\ref{PDCHDUId}; the seventh equality applies the mapping $\mathrm{ip}^{([1],X)@}$, introduced in Definition~\ref{DPTQIp}; the eighth equality follows from the fact that, by Definition~\ref{DIp}, $\mathrm{ip}^{(1,X)@}\circ\eta^{(1,X)}=\mathrm{ip}^{(1,X)}$; the ninth equality follows from the application of the mapping $\mathrm{CH}^{(1)\mathrm{m}}$ at a path class; the tenth equality follows from the fact that, by Proposition~\ref{PCHId}, $\mathrm{CH}^{(1)}\circ\mathrm{ip}^{(1,X)}=\eta^{(1,X)}$; finally, the last equality follows from the fact that, by Proposition~\ref{PDUEmb}, $\eta^{(2,1)\sharp}\circ\eta^{(1,X)}=\eta^{(2,X)}$.

This finishes the proof.
\end{proof}

Let us note that the second-order Curry-Howard mapping is also an $S$-sorted mapping from the underlying $S$-sorted set of the  partial $\Sigma^{\boldsymbol{\mathcal{A}}^{(2)}}$-algebra $\mathbf{Pth}_{\boldsymbol{\mathcal{A}}^{(2)}}$ to the underlying $S$-sorted set of the free $\Sigma^{\boldsymbol{\mathcal{A}}^{(2)}}$-algebra $\mathbf{T}_{\Sigma^{\boldsymbol{\mathcal{A}}^{(2)}}}(X)$. The question of whether $\mathrm{CH}^{(2)}$ is a $\Sigma^{\boldsymbol{\mathcal{A}}^{(2)}}$-homomorphism immediately raises. In this regard, we next prove that this mapping preserves, at least, the constants from the $S$-sorted sets $\mathcal{A}$ and $\mathcal{A}^{(2)}$.

\begin{figure}
\begin{center}
\begin{tikzpicture}
[ACliment/.style={-{To [angle'=45, length=5.75pt, width=4pt, round]}},scale=0.8]
\node[] (a) at (-4,-2) [] {$\mathcal{A}$};
\node[] (pth12) at (2,-2) [] {$\mathrm{Pth}_{\boldsymbol{\mathcal{A}}}$};
\node[] (pt12) at (2,-4) [] {$\mathrm{PT}_{\boldsymbol{\mathcal{A}}}$};
\node[] (pth1) at (8,-2) [] {$[\mathrm{Pth}_{\boldsymbol{\mathcal{A}}}]$};
\node[] (pt1) at (8,-4) [] {$[\mathrm{PT}_{\boldsymbol{\mathcal{A}}}]$};
\node[] (pth2) at (8,-6) [] {$\mathrm{Pth}_{\boldsymbol{\mathcal{A}}^{(2)}}$};
\node[] (t2) at (8,-8) [] {$\mathrm{T}_{\Sigma^{\boldsymbol{\mathcal{A}}^{(2)}}}(X)$};

\draw[ACliment]  (a) to node [above right]
{$\textstyle \mathrm{ech}^{(1,\mathcal{A})}$} (pth12);
\draw[ACliment, bend right=10]  (a) to node [above right]
{$\textstyle \eta^{(1,\mathcal{A})}$} (pt12);
\draw[ACliment, bend right=20]  (a) to node  [midway, fill=white, pos=0.85]
{$\textstyle \mathrm{ech}^{(2,\mathcal{A})}$} (pth2);
\draw[ACliment, bend right=30]  (a) to node  [below left, pos=0.7]
{$\textstyle \eta^{(2,\mathcal{A})}$} (t2);

\draw[ACliment, bend right=30]  (pt12) to node  [midway, fill=white]
{$\textstyle \eta^{(2,1)\sharp}$} (t2);

\draw[ACliment]  (pth12) to node [above]
{$\mathrm{pr}^{\mathrm{Ker}(\mathrm{CH}^{(1)})}$} (pth1);

\draw[ACliment]  (pt12) to node [below]
{$\mathrm{pr}^{\Theta^{[1]}}$} (pt1);

\draw[ACliment]  (pth1) to node [midway, fill=white]
{$\mathrm{CH}^{(1)\mathrm{m}}$} (pt12);

\draw[ACliment] 
($(pt1)+(0,-.35)$) to node [right] {$\textstyle \mathrm{ip}^{(2,[1])\sharp}$}  ($(pth2)+(0,.35)$);

\draw[ACliment] 
($(pth1)+(-.15,-.35)$) to node [left] {$\textstyle \mathrm{CH}^{[1]}$}  ($(pt1)+(-.15,.35)$);
\draw[ACliment] 
($(pt1)+(.15,+.35)$) to node [right] {$\textstyle \mathrm{ip}^{([1],X)@}$}  ($(pth1)+(.15,-.35)$);

\draw[ACliment] 
($(pth12)+(-.15,-.35)$) to node [left] {$\textstyle \mathrm{CH}^{(1)}$}  ($(pt12)+(-.15,.35)$);
\draw[ACliment] 
($(pt12)+(.15,+.35)$) to node [right] {$\textstyle \mathrm{ip}^{(1,X)@}$}  ($(pth12)+(.15,-.35)$);

\draw[ACliment] 
($(pth2)+(-.15,-.35)$) to node [left] {$\textstyle \mathrm{CH}^{(2)}$}  ($(t2)+(-.15,.35)$);

\end{tikzpicture}
\end{center}
\caption{Behaviour of $\mathrm{CH}^{(2)}$ relative to $\mathcal{A}$ at layers 1 \& 2.}
\label{FDCHA}
\end{figure}

\begin{proposition}\label{PDCHA} 
For the second-order Curry-Howard mapping we have that
\begin{itemize}
\item[(i)] $\mathrm{CH}^{(2)}\circ\mathrm{ech}^{(2,\mathcal{A})}=\eta^{(2,\mathcal{A})}.$
\end{itemize}
In particular, for every sort $s\in S$, and every rewrite rule $\mathfrak{p}\in\mathcal{A}_{s}$, we have that 
\begin{itemize}
\item[(ii)] $\mathrm{CH}^{(2)}_{s}(\mathfrak{p}^{\mathbf{Pth}_{\boldsymbol{\mathcal{A}}^{(2)}}})
=
\mathfrak{p}^{\mathbf{T}_{\Sigma^{\boldsymbol{\mathcal{A}}^{(2)}}}(X)}.
$
\end{itemize}
\end{proposition}
\begin{proof}
The reader is advised to consult the diagram presented in Figure~\ref{FDCHA}.

Let $s$ be a sort in $S$ and let $\mathfrak{p}$ be a rewrite rule in $\mathcal{A}_{s}$, then the following chain of equalities holds
\allowdisplaybreaks
\begin{align*}
\mathrm{CH}^{(2)}_{s}\left(
\mathrm{ech}^{(2,\mathcal{A})}_{s}\left(
\mathfrak{p}
\right)\right)
&=
\mathrm{CH}^{(2)}_{s}\left(
\mathrm{ip}^{(2,[1])\sharp}_{s}\left(
\left[\mathfrak{p}^{\mathbf{PT}_{\boldsymbol{\mathcal{A}}}}
\right]_{s}
\right)\right)
\tag{1}
\\&=
\eta^{(2,1)\sharp}_{s}\left(
\mathrm{CH}^{(1)\mathrm{m}}_{s}\left(
\mathrm{ip}^{([1],X)@}_{s}\left(
\left[\mathfrak{p}^{\mathbf{PT}_{\boldsymbol{\mathcal{A}}}}
\right]_{s}
\right)\right)\right)
\tag{2}
\\&=
\eta^{(2,1)\sharp}_{s}\left(
\mathrm{CH}^{(1)\mathrm{m}}_{s}\left(
\mathrm{ip}^{([1],X)@}_{s}\left(
\mathfrak{p}^{
[\mathbf{PT}_{\boldsymbol{\mathcal{A}}}]}
\right)\right)\right)
\tag{3}
\\&=
\eta^{(2,1)\sharp}_{s}\left(
\mathrm{CH}^{(1)\mathrm{m}}_{s}\left(
\mathfrak{p}^{[\mathbf{Pth}_{\boldsymbol{\mathcal{A}}}]}
\right)\right)
\tag{4}
\\&=
\eta^{(2,1)\sharp}_{s}\left(
\mathrm{CH}^{(1)\mathrm{m}}_{s}\left(
\left[\mathrm{ech}^{(1,\mathcal{A})}_{s}\left(
\mathfrak{p}
\right)\right]_{s}
\right)\right)
\tag{5}
\\&=
\eta^{(2,1)\sharp}_{s}\left(
\mathrm{CH}^{(1)}_{s}\left(
\mathrm{ech}^{(1,\mathcal{A})}_{s}\left(
\mathfrak{p}
\right)\right)\right)
\tag{6}
\\&=
\eta^{(2,1)\sharp}_{s}\left(
\eta^{(1,\mathcal{A})}_{s}\left(
\mathfrak{p}
\right)\right)
\tag{7}
\\&=
\eta^{(2,1)\sharp}_{s}\left(
\mathfrak{p}^{\mathbf{PT}_{\boldsymbol{\mathcal{A}}}}
\right)
\tag{8}
\\&=
\mathfrak{p}^{
\mathbf{T}_{\Sigma^{\boldsymbol{\mathcal{A}}^{(2)}}}(X)
}
\tag{9}
\\&=
\eta^{(2,\mathcal{A})}_{s}(\mathfrak{p}).
\tag{10}
\end{align*}

The first equality applies the definition of the mapping $\mathrm{ech}^{(2,\mathcal{A})}$ introduced in Definition~\ref{DDPth}; the second equality follows from the characterization of the value for the second-order Curry-Howard mapping of $(2,[1])$-identity second-order paths introduced in  Proposition~\ref{PDCHDUId}; the third equality recovers the interpretation of the constant symbol $\mathfrak{p}$ in the partial many-sorted $\Sigma^{\boldsymbol{\mathcal{A}}}$-algebra $[\mathbf{PT}_{\boldsymbol{\mathcal{A}}}]$ introduced in Proposition~\ref{PPTQCatAlg}; the fourth equality holds from the fact that, in virtue of Theorem~\ref{TIso}, $\mathrm{ip}^{([1],X)@}$ is a $\Sigma^{\boldsymbol{\mathcal{A}}}$-homomorphism from $[\mathbf{PT}_{\boldsymbol{\mathcal{A}}}]$ to $[\mathbf{Pth}_{\boldsymbol{\mathcal{A}}}]$; the fifth equality recovers the interpretation of the constant symbol $\mathfrak{p}$ in the partial many-sorted $\Sigma^{\boldsymbol{\mathcal{A}}}$-algebra $[\mathbf{Pth}_{\boldsymbol{\mathcal{A}}}]$ introduced in Proposition~\ref{PCHCatAlg}; the sixth equality applies the value of the mapping $\mathrm{CH}^{(1)\mathrm{m}}$ at a path class; the seventh equality follows from the fact that, according to Proposition~\ref{PCHA}, $\mathrm{CH}^{(1)}\circ\mathrm{ech}^{(1,\mathcal{A})}=\eta^{(1,\mathcal{A})}$; the eighth equality recovers the interpretation of the constant operation symbol $\mathfrak{p}$ in the partial many-sorted $\Sigma^{\boldsymbol{\mathcal{A}}}$-algebra $\mathbf{PT}_{\boldsymbol{\mathcal{A}}}$ introduced in Proposition~\ref{PPTCatAlg}; the ninth equality follows from the fact that, according to Proposition~\ref{PDUEmb}, $\eta^{(2,1)\sharp}$ is a $\Sigma^{\boldsymbol{\mathcal{A}}}$-homomorphism from $\mathbf{PT}_{\boldsymbol{\mathcal{A}}}$ to $\mathbf{T}_{\Sigma^{\boldsymbol{\mathcal{A}}^{(2)}}}(X)$; finally, the last equality recovers the value of the mapping $\eta^{(2,\mathcal{A})}$ at a rewrite rule as it is presented in Definition~\ref{DDEta}.

Let us note that, from the above equations, and taking into account Definitions~\ref{DDPth} and~\ref{DDEta}, it follows that, for every sort $s\in S$ and every rewrite rule $\mathfrak{p}$ in $\mathcal{A}_{s}$ we have that
$$
\mathrm{CH}^{(2)}_{s}\left(
\mathfrak{p}^{\mathbf{Pth}_{\boldsymbol{\mathcal{A}}^{(2)}}}
\right)
=
\mathfrak{p}^{\mathbf{T}_{\Sigma^{\boldsymbol{\mathcal{A}}^{(2)}}}(X)}.
$$

The statement follows.
\end{proof}

\begin{figure}
\begin{center}
\begin{tikzpicture}
[ACliment/.style={-{To [angle'=45, length=5.75pt, width=4pt, round]}
}, scale=0.8]
\node[] (xq) 		at 	(0,0) 	[] 	{$\mathcal{A}^{(2)}$};
\node[] (txq) 	at 	(6,0) 	[] 	{$\mathrm{Pth}_{\boldsymbol{\mathcal{A}}^{(2)}}$};
\node[] (txqc) 	at 	(6,-3) 	[] 	 {$\mathrm{T}_{\Sigma^{\boldsymbol{\mathcal{A}}^{(2)}}}(X)$};
\draw[ACliment]  (xq) 	to node [above]	
{$\mathrm{ech}^{(2,\mathcal{A}^{(2)})}$} (txq);
\draw[ACliment]  (txq) 	to node [right]	
  {$\mathrm{CH}^{(2)}$} (txqc);
\draw[ACliment, bend right=10]  (xq) 	to node [below left]	
{$\eta^{(2,\mathcal{A}^{(2)})}$} (txqc);
\end{tikzpicture}
\end{center}
\caption{Behaviour of $\mathrm{CH}^{(2)}$ relative to $\mathcal{A}^{(2)}$ at layer 2.}
\label{FDCHDA}
\end{figure}

\begin{proposition}\label{PDCHDA} 
For the second-order Curry-Howard mapping we have that
\begin{itemize}
\item[(i)] $\mathrm{CH}^{(2)}\circ\mathrm{ech}^{(2,\mathcal{A}^{(2)})}=\eta^{(2,\mathcal{A}^{(2)})}.$
\end{itemize}
In particular, for every sort $s\in S$, and every second-order rewrite rule $\mathfrak{p}^{(2)}\in\mathcal{A}^{(2)}_{s}$, we have that 
\begin{itemize}
\item[(ii)] $\mathrm{CH}^{(2)}_{s}(\mathfrak{p}^{(2)
\mathbf{Pth}_{\boldsymbol{\mathcal{A}}^{(2)}}
})
=
\mathfrak{p}^{(2)\mathbf{T}_{\Sigma^{\boldsymbol{\mathcal{A}}^{(2)}}}(X)}.
$
\end{itemize}
\end{proposition}
\begin{proof}
The reader is advised to consult the diagram presented in Figure~\ref{FDCHDA}.

Let $s$ be a sort in $S$ and $\mathfrak{p}^{(2)}$ a second-order rewrite rule in $\mathcal{A}^{(2)}_{s}$, then the following chain of equalities holds
\begin{align*}
\mathrm{CH}^{(2)}_{s}\left(
\mathrm{ech}^{(2,\mathcal{A})^{(2)}}_{s}
\left(
\mathfrak{p}^{(2)}
\right)
\right)
&=
\mathrm{CH}^{(2)}_{s}\left(
\mathfrak{p}^{(2)\mathbf{Pth}_{\boldsymbol{\mathcal{A}}^{(2)}}}
\right)
\tag{1}
\\&=
\mathfrak{p}^{(2)\mathbf{T}_{\Sigma^{\boldsymbol{\mathcal{A}}^{(2)}}}(X)}
\tag{2}
\\&=
\eta^{(2,\mathcal{A}^{(2)})}_{s}(\mathfrak{p}^{(2)}).
\tag{3}
\end{align*}

The first equality recovers the interpretation of the constant symbol $\mathfrak{p}^{(2)}$ in the partial $\Sigma^{\boldsymbol{\mathcal{A}}^{(2)}}$-algebra $\mathbf{Pth}_{\boldsymbol{\mathcal{A}}^{(2)}}$ introduced in Proposition~\ref{PDPthDCatAlg}; the second equality follows by Definition~\ref{DDCH}, since the second-order Curry-Howard mapping assigns to a second-order echelon the unique second-order rewrite rule appearing in it; finally, the last equality recovers the interpretation of the mapping $\eta^{(2,\mathcal{A}^{(2)})}$ at a second-order rewrite rule $\mathfrak{p}^{(2)}$ as it is presented in Definition~\ref{DDEta}.

This finishes the proof.
\end{proof}

However the preservation of the constants from $\mathcal{A}$ and $\mathcal{A}^{(2)}$ is not enough to infer  that the second-order Curry-Howard mapping is a  $\Sigma^{\boldsymbol{\mathcal{A}}^{(2)}}$-homomorphism. This is due, essentially, to the fact that two things occur here. On one hand the many-sorted free $\Sigma^{\boldsymbol{\mathcal{A}}^{(2)}}$-algebra $\mathbf{T}_{\Sigma^{\boldsymbol{\mathcal{A}}^{(2)}}}(X)$ carries the same problems as the many-sorted free $\Sigma^{\boldsymbol{\mathcal{A}}}$-algebra $\mathrm{T}_{\Sigma^{\boldsymbol{\mathcal{A}}}}(X)$ in the first part. Note that, according to Proposition~\ref{PDCHDUId} $\mathrm{CH}^{(2)}$ behaves essentially equal to $\mathrm{CH}^{(1)}$ on $(2,[1])$-identity second-order paths, and $\mathrm{CH}^{(1)}$ was not a $\Sigma^{\boldsymbol{\mathcal{A}}}$-homomorphism as it was shown in Proposition~\ref{PCHNotHomCat}. This problem was solved in the first part of this work by constructing a $\Sigma^{\boldsymbol{\mathcal{A}}}$-congruence on $\mathbf{T}_{\Sigma^{\boldsymbol{\mathcal{A}}}}(X)$, i.e., $\Theta^{[1]}$, that solved these problems. We will see later that this approach can also be used here.

On the other hand, this second iteration will have the same problems on non-$(2,[1])$-identity second-order paths, since  several terms can denote the same second-order path. This problem will be later addressed  following a similar criterion.

We next present two counterexamples along the lines described above showing that $\mathrm{CH}^{(2)}$ is neither a $\Sigma^{\boldsymbol{\mathcal{A}}}$-homomorphism nor a $\Sigma^{\boldsymbol{\mathcal{A}}^{(2)}}$-homomorphism.

\begin{restatable}{proposition}{PDCHNotHomCat}
\label{PDCHNotHomCat}
The second-order Curry-Howard mapping $\mathrm{CH}^{(2)}$ is not a $\Sigma^{\boldsymbol{\mathcal{A}}}$-homomorphism from the partial $\Sigma^{\boldsymbol{\mathcal{A}}}$-algebra $\mathbf{Pth}^{(1,2)}_{\boldsymbol{\mathcal{A}}^{(2)}}$ to the $\Sigma^{\boldsymbol{\mathcal{A}}}$-algebra $\mathbf{T}^{(1,2)}_{\Sigma^{\boldsymbol{\mathcal{A}}^{(2)}}}(X)$.
\end{restatable}
\begin{proof}
For a non-empty set of sorts $S$, consider a non-empty $S$-sorted set $X$ and let $\Sigma$ be the empty $S$-sorted signature. Let $s$ be a sort in $S$ and $x$ a variable in $X_{s}$. Let us consider the $(2,0)$-identity second-order path on $x$, that is, consider the second-order path $\mathfrak{P}^{(2)}=\mathrm{ip}^{(2,X)}_{s}(x)$.

The following chain of equalities holds.
\allowdisplaybreaks
\begin{align*}
\mathrm{sc}^{(0,2)}_{s}\left(
\mathfrak{P}^{(2)}\right)&=
\mathrm{sc}^{(0,2)}_{s}\left(
\mathrm{ip}^{(2,[1])\sharp}_{s}\left(
\left[\eta^{(1,X)}_{s}(x)
\right]_{s}
\right)\right)
\tag{1}
\\&=
\mathrm{sc}^{(0,[1])}_{s}\left(
\mathrm{ip}^{([1],X)@}_{s}\left(
\mathrm{sc}^{([1],2)}_{s}\left(
\mathrm{ip}^{(2,[1])\sharp}_{s}\left(
\left[\eta^{(1,X)}_{s}(x)
\right]_{s}
\right)\right)\right)\right)
\tag{2}
\\&=
\mathrm{sc}^{(0,[1])}_{s}\left(
\mathrm{ip}^{([1],X)@}_{s}\left(
\left[\eta^{(1,X)}_{s}(x)
\right]_{s}
\right)\right)
\tag{3}
\\&=
\mathrm{sc}^{(0,[1])}_{s}\left(
\left[
\mathrm{ip}^{(1,X)@}_{s}\left(
\eta^{(1,X)}_{s}(x)
\right)
\right]_{s}
\right)
\tag{4}
\\&=
\mathrm{sc}^{(0,[1])}_{s}\left(
\left[
\mathrm{ip}^{(1,X)}_{s}(x)
\right]_{s}
\right)
\tag{5}
\\&=
\mathrm{sc}^{(0,1)}_{s}\left(
\mathrm{ip}^{(1,X)}_{s}\left(
x
\right)\right)
\tag{6}
\\&=
\eta^{(0,X)}_{s}(x).
\tag{7}
\end{align*}

The first equality follows from the characterization of $\mathfrak{P}^{(2)}$ presented in Definition~\ref{DDPth}; the second equality unravels the mapping $\mathrm{sc}^{(0,2)}$ introduced in Definition~\ref{DDScTgZ}; the third equality follows from Proposition~\ref{PDBasicEq}; the fourth equality unravels the mapping $\mathrm{ip}^{([1],X)@}$ introduced in Definition~\ref{DPTQIp}; the fifth equality follows from Definition~\ref{DIp}; the sixth equality applies the mapping $\mathrm{sc}^{(0,[1])}$ at a path class as it is presented in Definition~\ref{DCHUZ}; finally, the last equality follows from Proposition~\ref{PBasicEqX}.

By a similar sequence of arguments we can also infer that
$$
\mathrm{tg}^{(0,2)}_{s}\left(
\mathfrak{P}^{(2)}
\right)=\eta^{(0,X)}_{s}(x).
$$

It follows from Claim~\ref{CDPthCatAlgCompZ} that the $0$-composition of $\mathfrak{P}^{(2)}$ with itself is well-defined. Since $\mathfrak{P}^{(2)}$ is a $(2,[1])$-identity second-order path, we have that its $0$-composite is a $(2,[1])$-identity second-order path as well.
In fact, from Claim~\ref{CDPthCatAlgCompZ} and Proposition~\ref{PDUIp}, we have that
$$
\mathfrak{P}^{(2)}\circ^{0\mathbf{Pth}_{\boldsymbol{\mathcal{A}}^{(2)}}}_{s}\mathfrak{P}^{(2)}
=
\mathrm{ip}^{(2,[1])\sharp}_{s}\left(
\left[
\eta^{(1,X)}_{s}(x)
\circ_{s}^{0\mathbf{PT}_{\boldsymbol{\mathcal{A}}}}
\eta^{(1,X)}_{s}(x)
\right]_{s}
\right).
$$

Note that the following chain of equalities holds
\begin{flushleft}
$\mathrm{ip}^{(1,X)@}_{s}\left(
\eta^{(1,X)}_{s}(x)
\circ^{0\mathbf{PT}_{\boldsymbol{\mathcal{A}}}}_{s}
\eta^{(1,X)}_{s}(x)
\right)$
\allowdisplaybreaks
\begin{align*}
&=
\mathrm{ip}^{(1,X)@}_{s}\left(
\eta^{(1,X)}_{s}(x)
\right)
\circ^{0\mathbf{Pth}_{\boldsymbol{\mathcal{A}}}}_{s}
\mathrm{ip}^{(1,X)@}_{s}\left(
\eta^{(1,X)}_{s}(x)
\right)
\tag{1}
\\&=
\mathrm{ip}^{(1,X)}_{s}(x)
\circ_{s}^{0\mathbf{Pth}_{\boldsymbol{\mathcal{A}}}}
\mathrm{ip}^{(1,X)}_{s}(x)
\tag{2}
\\&=
\mathrm{ip}^{(1,X)}_{s}(x).
\tag{3}
\end{align*}
\end{flushleft}

The first equality follows from the fact that, by Definition~\ref{DIp}, $\mathrm{ip}^{(1,X)@}$ is a $\Sigma^{\boldsymbol{\mathcal{A}}}$-homomorphism; the second equality follows from the Universal Property of $\mathrm{ip}^{(1,X)@}$ presented in Definition~\ref{DIp}; finally, the last equality follows from the fact that $(1,0)$-identity paths are idempotent for the $0$-composition.

Note also that, according to Proposition~\ref{PCHId}, we have that 
$$
\mathrm{CH}^{(1)}_{s}\left(
\mathrm{ip}^{(1,X)}_{s}(
x
)\right)
=
\eta^{(1,X)}_{s}(x).
$$

Thus, by Lemma~\ref{LWCong}, we have that 
$$
\left[
\eta^{(1,X)}_{s}(x)
\circ_{s}^{0\mathbf{PT}_{\boldsymbol{\mathcal{A}}}}
\eta^{(1,X)}_{s}(x)
\right]_{s}
=
\left[\eta^{(1,X)}_{s}(x)
\right]_{s}.
$$

Thus, $\mathfrak{P}^{(2)}
\circ_{s}^{0\mathbf{Pth}_{\boldsymbol{\mathcal{A}}^{(2)}}}
\mathfrak{P}^{(2)}=\mathfrak{P}^{(2)}
$, that is $\mathfrak{P}^{(2)}$ is idempotent for the $0$-composition operation on second-order paths.

Thus, the following chain of equalities holds
\allowdisplaybreaks
\begin{align*}
\mathrm{CH}^{(2)}_{s}\left(
\mathfrak{P}^{(2)}
\circ_{s}^{0\mathbf{Pth}_{\boldsymbol{\mathcal{A}}^{(2)}}}
\mathfrak{P}^{(2)}
\right)
&=
\mathrm{CH}^{(2)}_{s}\left(
\mathfrak{P}^{(2)}
\right)
\tag{1}
\\&=
\mathrm{CH}^{(2)}_{s}\left(
\mathrm{ip}^{(2,[1])\sharp}_{s}\left(
\left[\eta^{(1,X)}_{s}(x)
\right]_{s}
\right)\right)
\tag{2}
\\&=
\eta^{(2,1)\sharp}_{s}\left(
\mathrm{CH}^{(1)\mathrm{m}}_{s}\left(
\mathrm{ip}^{([1],X)@}_{s}\left(
\left[\eta^{(1,X)}_{s}(x)
\right]_{s}
\right)\right)\right)
\tag{3}
\\&=
\eta^{(2,1)\sharp}_{s}\left(
\mathrm{CH}^{(1)\mathrm{m}}_{s}\left(
\left[\mathrm{ip}^{(1,X)@}_{s}\left(
\eta^{(1,X)}_{s}\left(
x\right)\right)
\right]_{s}
\right)\right)
\tag{4}
\\&=
\eta^{(2,1)\sharp}_{s}\left(
\mathrm{CH}^{(1)\mathrm{m}}_{s}\left(
\left[\mathrm{ip}^{(1,X)}_{s}(
x)
\right]_{s}
\right)\right)
\tag{5}
\\&=
\eta^{(2,1)\sharp}_{s}\left(
\mathrm{CH}^{(1)}_{s}\left(
\mathrm{ip}^{(1,X)}_{s}(
x
)\right)\right)
\tag{6}
\\&=
\eta^{(2,1)\sharp}_{s}\left(
\eta^{(1,X)}_{s}(
x
)\right)
\tag{7}
\\&=
\eta^{(2,X)}_{s}(x).
\tag{8}
\end{align*}

The first equality follows from the previous discussion; the second equality follows the characterization of $\mathfrak{P}^{(2)}$ introduced in Definition~\ref{DDPth}; the third equality follows from Proposition~\ref{PDCHDUId}; the fourth equality applies the mapping $\mathrm{ip}^{([1],X)@}$ as it is introduced in Definition~\ref{DPTQIp}; the fifth equality follows from the Universal Property of $\mathrm{ip}^{(1,X)@}$ as it is presented in Definition~\ref{DIp}; the sixth equality applies the monomorphic Curry-Howard mapping at a path class; the seventh equality follows from Proposition~\ref{PCHId}; finally, the last equality follows from Proposition~\ref{PDUEmb}.

If we assume, towards a contradiction, that $\mathrm{CH}^{(2)}$ is a $\Sigma^{\boldsymbol{\mathcal{A}}}$-homomorphism. Then we would obtain the following  chain of equalities
\allowdisplaybreaks
\begin{align*}
\mathrm{CH}^{(2)}_{s}\left(
\mathfrak{P}^{(2)}
\circ_{s}^{0\mathbf{Pth}_{\boldsymbol{\mathcal{A}}^{(2)}}}
\mathfrak{P}^{(2)}
\right)
&=
\mathrm{CH}^{(2)}_{s}\left(
\mathfrak{P}^{(2)}
\right)
\circ_{s}^{0\mathbf{T}_{\Sigma^{\boldsymbol{\mathcal{A}}^{(2)}}}(X)}
\mathrm{CH}^{(2)}_{s}\left(
\mathfrak{P}^{(2)}
\right)
\tag{1}
\\&=
\eta^{(2,X)}_{s}(x)
\circ_{s}^{0\mathbf{T}_{\Sigma^{\boldsymbol{\mathcal{A}}^{(2)}}}(X)}
\eta^{(2,X)}_{s}(x).
\tag{2}
\end{align*}

The first equality would follow from the assumption of $\mathrm{CH}^{(2)}$ being a $\Sigma^{\boldsymbol{\mathcal{A}}}$-homomorphism; whilst the last equality would follow from the previous discussion on the value of the second-order Curry-Howard mapping at $\mathfrak{P}^{(2)}$.

However the terms $\eta^{(2,X)}_{s}(x)
\circ_{s}^{0\mathbf{T}_{\Sigma^{\boldsymbol{\mathcal{A}}^{(2)}}}(X)}
\eta^{(2,X)}_{s}(x)$  and $\eta^{(2,X)}_{s}(x)$ are different in $\mathbf{T}_{\Sigma^{\boldsymbol{\mathcal{A}}^{(2)}}}(X)_{s}$.
\end{proof}

As an immediate consequence of the previous result, we have that $\mathrm{CH}^{(2)}$ is not a $\Sigma^{\boldsymbol{\mathcal{A}}^{(2)}}$-homomorphism from the partial $\Sigma^{\boldsymbol{\mathcal{A}}^{(2)}}$-algebra $\mathbf{Pth}_{\boldsymbol{\mathcal{A}}^{(2)}}$ to the free $\Sigma^{\boldsymbol{\mathcal{A}}^{(2)}}$-algebra $\mathbf{T}_{\Sigma^{\boldsymbol{\mathcal{A}}^{(2)}}}(X)$. However, in the next proposition, we will introduce a second-order path much more aligned with the counterexample  introduced in Proposition~\ref{PCHNotHomCat}.

\begin{restatable}{proposition}{PDCHNotHomDCat}
\label{PDCHNotHomDCat}
The second-order Curry-Howard mapping $\mathrm{CH}^{(2)}$ is not a $\Sigma^{\boldsymbol{\mathcal{A}}^{(2)}}$-homomorphism from the partial $\Sigma^{\boldsymbol{\mathcal{A}}^{(2)}}$-algebra $\mathbf{Pth}_{\boldsymbol{\mathcal{A}}^{(2)}}$ to the free $\Sigma^{\boldsymbol{\mathcal{A}}^{(2)}}$-algebra $\mathbf{T}_{\Sigma^{\boldsymbol{\mathcal{A}}^{(2)}}}(X)$.
\end{restatable}
\begin{proof}
Let $s$ be a sort in $S$ and $[P]_{s}$ a path term class in $[\mathrm{PT}_{\boldsymbol{\mathcal{A}}}]_{s}$. Let us consider the $(2,[1])$-identity path on $[P]_{s}$, i.e., $\mathrm{ip}^{(2,[1])\sharp}_{s}([P]_{s})$. In virtue of Proposition~\ref{PDCHDUId}, we have that
$$
\mathrm{CH}^{(2)}_{s}\left(
\mathrm{ip}^{(2,[1])\sharp}_{s}\left(
\left[
P
\right]_{s}
\right)\right)=
\eta^{(2,1)\sharp}_{s}\left(
\mathrm{CH}^{(1)\mathrm{m}}_{s}\left(
\mathrm{ip}^{([1],X)@}_{s}\left(
\left[P
\right]_{s}
\right)\right)\right)
.
$$

However, the $(2,[1])$-identity path on $[P]_{s}$ can be $1$-composed with itself and we have 
$$
\mathrm{ip}^{(2,[1])\sharp}_{s}\left(
\left[
P
\right]_{s}
\right)
\circ_{s}^{1\mathbf{Pth}_{\boldsymbol{\mathcal{A}}^{(2)}}}
\mathrm{ip}^{(2,[1])\sharp}_{s}\left(
\left[P
\right]_{s}
\right)
=\mathrm{ip}^{(2,[1])\sharp}_{s}\left(
\left[P
\right]_{s}\right).
$$

Let us assume, towards a contradiction, that $\mathrm{CH}^{(2)}$ is a $\Sigma^{\boldsymbol{\mathcal{A}}^{(2)}}$-homomorphism. Then we have the following chain of equalities
\begin{flushleft}
$\eta^{(2,1)\sharp}_{s}\left(
\mathrm{CH}^{(1)\mathrm{m}}_{s}\left(
\mathrm{ip}^{([1],X)@}_{s}\left(
\left[
P
\right]_{s}
\right)\right)\right)$
\begin{align*}
&=\mathrm{CH}^{(2)}_{s}
\left(
\mathrm{ip}^{(2,[1])\sharp}_{s}
\left(
\left[
P
\right]_{s}
\right)\right)
\tag{1}
\\&=\mathrm{CH}^{(2)}_{s}\left(
\mathrm{ip}^{(2,[1])\sharp}_{s}\left(
\left[P
\right]_{s}
\right)
\circ_{s}^{1\mathbf{Pth}_{\boldsymbol{\mathcal{A}}^{(2)}}}
\mathrm{ip}^{(2,[1])\sharp}_{s}
\left(
\left[P
\right]_{s}
\right)\right)
\tag{2}
\\&=
\mathrm{CH}^{(2)}_{s}
\left(
\mathrm{ip}^{(2,[1])\sharp}_{s}
\left(
\left[
P
\right]_{s}
\right)\right)
\circ_{s}^{1\mathbf{T}_{\Sigma^{\boldsymbol{\mathcal{A}}^{(2)}}}(X)}
\mathrm{CH}^{(2)}_{s}\left(
\mathrm{ip}^{(2,[1])\sharp}_{s}
\left(\left[
P
\right]_{s}
\right)\right)
\tag{3}
\\&=
\eta^{(2,1)\sharp}_{s}\left(
\mathrm{CH}^{(1)\mathrm{m}}_{s}\left(
\mathrm{ip}^{([1],X)@}_{s}\left(
\left[
P
\right]_{s}
\right)\right)\right)
\circ_{s}^{1\mathbf{T}_{\Sigma^{\boldsymbol{\mathcal{A}}^{(2)}}}(X)}
\\&
\qquad\qquad\qquad\qquad\qquad\qquad
\qquad
\eta^{(2,1)\sharp}_{s}\left(
\mathrm{CH}^{(1)\mathrm{m}}_{s}\left(
\mathrm{ip}^{([1],X)@}_{s}\left(
\left[
P
\right]_{s}
\right)\right)\right).
\tag{4}
\end{align*}
\end{flushleft}

The first equality follows from Proposition~\ref{PDCHDUId}; the second equality follows from the fact that $\mathrm{ip}^{(2,[1])\sharp}_{s}([P]_{s})$ is idempotent for the $1$-composition; the third equality follows since we are assuming that $\mathrm{CH}^{(2)}$ is a $\Sigma^{\boldsymbol{\mathcal{A}}^{(2)}}$-homomorphism; finally the last equality follows from Proposition~\ref{PDCHDUId}.

However, the terms $\eta^{(2,1)\sharp}_{s}(
\mathrm{CH}^{(1)\mathrm{m}}_{s}(
\mathrm{ip}^{([1],X)@}_{s}(
[P]_{s}
)))$ and
\[
\eta^{(2,1)\sharp}_{s}\left(
\mathrm{CH}^{(1)\mathrm{m}}_{s}\left(
\mathrm{ip}^{([1],X)@}_{s}\left(
[P]_{s}
\right)\right)\right)
\circ_{s}^{1\mathbf{T}_{\Sigma^{\boldsymbol{\mathcal{A}}^{(2)}}}(X)}
\eta^{(2,1)\sharp}_{s}\left(
\mathrm{CH}^{(1)\mathrm{m}}_{s}\left(
\mathrm{ip}^{([1],X)@}_{s}\left(
[P]_{s}
\right)\right)\right)
\]
are different  in $\mathrm{T}_{\Sigma^{\boldsymbol{\mathcal{A}}^{(2)}}}(X)_{s}$.
\end{proof}

Nevertheless, there are some specific configurations on  $1$-compositions of second-order paths for which the second-order Curry-Howard mapping behaves like a  $\Sigma^{\boldsymbol{\mathcal{A}}^{(2)}}$-homomorphism. Since later on we will make use of these cases, we introduce the following result.

\begin{proposition}\label{PDCHEchHom}
Let $s$ be a sort in $S$ and let $\mathfrak{P}^{(2)}$ be a second-order path in  $\mathrm{Pth}_{\boldsymbol{\mathcal{A}}^{(2)},s}$ of length strictly greater than one containing at least one echelon. Let $i\in\bb{\mathfrak{P}^{(2)}}$ be the first index for which $\mathfrak{P}^{(2),i,i}$ is an echelon. Then
\begin{itemize}
\item[(i)] if $i=0$, that is, if the first echelon appears in the first step, then
\begin{multline*}
\mathrm{CH}^{(2)}_{s}
\left(
\mathfrak{P}^{(2),1,\bb{\mathfrak{P}^{(2)}}-1}
\circ_{s}^{1\mathbf{Pth}_{\boldsymbol{\mathcal{A}}^{(2)}}}
\mathfrak{P}^{(2),0,0}
\right)
\\=
\mathrm{CH}^{(2)}_{s}
\left(
\mathfrak{P}^{(2),1,\bb{\mathfrak{P}^{(2)}}-1}
\right)
\circ_{s}^{1\mathbf{T}_{\Sigma^{\boldsymbol{\mathcal{A}}^{(2)}}}(X)}
\mathrm{CH}^{(2)}_{s}
\left(
\mathfrak{P}^{(2),0,0}
\right);
\end{multline*}
\item[(ii)] if $i=\bb{\mathfrak{P}^{(2)}}$, that is, if the first echelon appears in the last step, then
\begin{multline*}
\mathrm{CH}^{(2)}_{s}
\left(
\mathfrak{P}^{(2),\bb{\mathfrak{P}^{(2)}}-1,\bb{\mathfrak{P}^{(2)}}-1}
\circ_{s}^{1\mathbf{Pth}_{\boldsymbol{\mathcal{A}}^{(2)}}}
\mathfrak{P}^{(2),0,\bb{\mathfrak{P}^{(2)}}-2}
\right)
\\=
\mathrm{CH}^{(2)}_{s}
\left(
\mathfrak{P}^{(2),\bb{\mathfrak{P}^{(2)}}-1,\bb{\mathfrak{P}^{(2)}}-1}
\right)
\circ_{s}^{1\mathbf{T}_{\Sigma^{\boldsymbol{\mathcal{A}}^{(2)}}}(X)}
\mathrm{CH}^{(2)}_{s}
\left(
\mathfrak{P}^{(2),0,\bb{\mathfrak{P}^{(2)}}-2}
\right).
\end{multline*}
\end{itemize}
\end{proposition}
\begin{proof}
These cases were already considered in Definition~\ref{DDCH}.
\end{proof}

The last case to consider is that of echelonless head-constant coherent second-order paths. 

\begin{proposition}\label{PDCHCat} Let $s$ be a sort in $S$ and let $\mathfrak{P}^{(2)}$ be an echelonless coherent head-constant second-order path in $\mathrm{Pth}_{\boldsymbol{\mathcal{A}}^{(2)},s}$. Let $\mathbf{s}\in S^{\star}-\{\lambda\}$ and $\sigma$ and operation symbol in $\Sigma_{\mathbf{s},s}$ be the unique word and operation symbol, respectively, associated with $\mathfrak{P}^{(2)}$. Let $(\mathfrak{P}^{(2)}_{j})_{j\in\bb{\mathbf{s}}}$ in $\mathrm{Pth}_{\boldsymbol{\mathcal{A}}^{(2)},\mathbf{s}}$ be the family of second-order path terms that we can extract from $\mathfrak{P}^{(2)}$ in virtue of Lemma~\ref{LDPthExtract}. The following equalities hold
\[
\mathrm{CH}^{(2)}_{s}\left(\mathfrak{P}^{(2)}\right)
=
\mathrm{CH}^{(2)}_{s}\left(\sigma^{\mathbf{Pth}_{\boldsymbol{\mathcal{A}}^{(2)}}}\left(\left(
\mathfrak{P}^{(2)}_{j}\right)_{j\in\bb{\mathbf{s}}}\right)\right).
\]
\end{proposition}
\begin{proof}
Following Definition~\ref{DDCH}, the value of the second-order Curry-Howard mapping at $\mathfrak{P}^{(2)}$ is given by
\[
\mathrm{CH}^{(2)}_{s}\left(\mathfrak{P}^{(2)}\right)
=\sigma^{\mathbf{T}_{\Sigma^{\boldsymbol{\mathcal{A}}^{(2)}}}(X)}\left(
\left(
\mathrm{CH}^{(2)}_{s_{j}}\left(
\mathfrak{P}^{(2)}_{j}
\right)
\right)_{j\in\bb{\mathbf{s}}}
\right).
\]

Note that, by assumption, $\mathfrak{P}^{(2)}$, being echelonless, is not a $(2,[1])$-identity second-order path. Thus, the family of second-order paths we can extract from $\mathfrak{P}^{(2)}$, i.e., $(\mathfrak{P}^{(2)}_{j})_{j\in\bb{\mathbf{s}}}$, is not a family of $(2,[1])$-identity second-order paths. Thus, according to Corollary~\ref{CDPthWB}, the second-order path $\sigma^{\mathbf{Pth}_{\boldsymbol{\mathcal{A}}^{(2)}}}((
\mathfrak{P}^{(2)}_{j})_{j\in\bb{\mathbf{s}}})$ is a coherent head-constant echelonless second-order path. Moreover, according to Proposition~\ref{PDRecov}, the family of second-order paths that we can extract from $\sigma^{\mathbf{Pth}_{\boldsymbol{\mathcal{A}}^{(2)}}}((
\mathfrak{P}^{(2)}_{j})_{j\in\bb{\mathbf{s}}})$ in virtue of Lemma~\ref{LDPthExtract} is, precisely, $(\mathfrak{P}^{(2)}_{j})_{j\in\bb{\mathbf{s}}}$.

Following Definition~\ref{DDCH}, the value of the second-order Curry-Howard mapping at $\sigma^{\mathbf{Pth}_{\boldsymbol{\mathcal{A}}^{(2)}}}((
\mathfrak{P}^{(2)}_{j})_{j\in\bb{\mathbf{s}}})$ is given by
\[
\mathrm{CH}^{(2)}_{s}\left(\sigma^{\mathbf{Pth}_{\boldsymbol{\mathcal{A}}^{(2)}}}\left(\left(
\mathfrak{P}^{(2)}_{j}\right)_{j\in\bb{\mathbf{s}}}\right)\right)
=\sigma^{\mathbf{T}_{\Sigma^{\boldsymbol{\mathcal{A}}^{(2)}}}(X)}\left(
\left(
\mathrm{CH}^{(2)}_{s_{j}}\left(
\mathfrak{P}^{(2)}_{j}
\right)
\right)_{j\in\bb{\mathbf{s}}}
\right).
\]

This completes the proof.
\end{proof}

As it will be of interest later, we give an end to this section proving that $\coprod \mathrm{CH}^{(2)}$ is a monotone mapping from the partially ordered set $(\coprod\mathrm{Pth}_{\boldsymbol{\mathcal{A}}^{(2)}}, \leq_{\mathbf{Pth}_{\boldsymbol{\mathcal{A}}^{(2)}}})$ to the partially ordered set $(\coprod\mathrm{T}_{\Sigma^{\boldsymbol{\mathcal{A}}^{(2)}}}(X), \leq_{\mathbf{T}_{\Sigma^{\boldsymbol{\mathcal{A}}^{(2)}}}(X)})$. This will ultimately be crucial, as more complex paths have more complex representations.

\begin{restatable}{proposition}{PDCHMono}
\label{PDCHMono} The mapping $\coprod\mathrm{CH}^{(2)}$ from $\coprod\mathrm{Pth}_{\boldsymbol{\mathcal{A}}^{(2)}}$ to $\coprod\mathrm{T}_{\Sigma^{\boldsymbol{\mathcal{A}}^{(2)}}}$ is order-preserving 
$$
\textstyle
\coprod\mathrm{CH}^{(2)}\colon
\left(
\coprod
\mathrm{Pth}_{\boldsymbol{\mathcal{A}}^{(2)}}, \leq_{
\mathbf{Pth}_{\boldsymbol{\mathcal{A}}^{(2)}}}
\right)
\mor
\left(\coprod\mathrm{T}_{\Sigma^{\boldsymbol{\mathcal{A}}^{(2)}}}(X), \leq_{\mathbf{T}_{\Sigma^{\boldsymbol{\mathcal{A}}^{(2)}}}(X)}
\right),
$$
that is, given pairs $(\mathfrak{Q}^{(2)},t)$, $(\mathfrak{P}^{(2)},s)$ in $\coprod\mathrm{Pth}_{\boldsymbol{\mathcal{A}}^{(2)}}$, if $(\mathfrak{Q}^{(2)},t)
\leq_{
\mathbf{Pth}_{\boldsymbol{\mathcal{A}}^{(2)}}}
(\mathfrak{P}^{(2)},s)
$, then $\mathrm{CH}^{(2)}_{t}(\mathfrak{Q}^{(2)})$ is a subterm of type $t$ of the term $\mathrm{CH}^{(2)}_{s}(\mathfrak{P}^{(2)})$, i.e., 
$$
\mathrm{CH}^{(2)}_{t}
\left(
\mathfrak{Q}^{(2)}
\right)
\in\mathrm{Subt}_{\Sigma^{\boldsymbol{\mathcal{A}}^{(2)}}}
\left(
\mathrm{CH}^{(2)}_{s}
\left(\mathfrak{P}^{(2)}
\right)\right)_{t}.
$$
\end{restatable}
\begin{proof}
Let us recall from Remark~\ref{RDOrd} that $(\mathfrak{Q}^{(2)},t)
\leq_{
\mathbf{Pth}_{\boldsymbol{\mathcal{A}}^{(2)}}}
(\mathfrak{P}^{(2)},s)$ if, and only if, $s=t$ and $\mathfrak{Q}^{(2)}=\mathfrak{P}^{(2)}$ or there exists a natural number $m\in\mathbb{N}-\{0\}$, a word $\mathbf{w}\in S^{\star}$ of length $\bb{\mathbf{w}}=m+1$, and a family of second-order paths $(\mathfrak{R}^{(2)}_{k})_{k\in\bb{\mathbf{w}}}$ in $\mathrm{Pth}_{\boldsymbol{\mathcal{A}}^{(2)},\mathbf{w}}$ such that $w_{0}=t$, $\mathfrak{R}^{(2)}_{0}=\mathfrak{Q}^{(2)}$, $w_{m}=s$, $\mathfrak{R}^{(2)}_{m}=\mathfrak{P}^{(2)}$ and, for every $k\in m$, $(\mathfrak{R}^{(2)}_{k},w_{k})\prec_{\mathbf{Pth}_{\boldsymbol{\mathcal{A}}^{(2)}}}
(\mathfrak{R}^{(2)}_{k+1},w_{k+1})
$.

The proposition holds trivially in case $s=t$ and $\mathfrak{Q}^{(2)}=\mathfrak{P}^{(2)}$. Therefore, it remains to prove the other case. We will prove it by induction on $m\in\mathbb{N}-\{0\}$.

\textsf{Base step of the induction.}

For $m=1$ we have that $(\mathfrak{Q}^{(2)},t)\prec_{\mathbf{Pth}_{\boldsymbol{\mathcal{A}}^{(2)}}} (\mathfrak{P}^{(2)},s)$, hence, according to Definition~\ref{DDOrd}, we are in one of the following situations
\begin{enumerate}
\item $\mathfrak{P}^{(2)}$ is a $(2,[1])$-identity second-order path; or
\item $\mathfrak{P}^{(2)}$ is a second-order path of length strictly greater than one containing at least one   second-order echelon; or
\item $\mathfrak{P}^{(2)}$ is an echelonless second-order path.
\end{enumerate}

If~(1), then $\mathfrak{P}^{(2)}$ and $\mathfrak{Q}^{(2)}$ are $(2,[1])$-identity second-order paths of the form
\begin{align*}
\mathfrak{P}^{(2)}&=
\mathrm{ip}^{(2,[1])\sharp}_{s}\left(
[P]_{s}
\right),
&
\mathfrak{Q}^{(2)}&=
\mathrm{ip}^{(2,[1])\sharp}_{t}\left(
[Q]_{t}
\right),
\end{align*}
for some path term classes $[P]_{s}\in [\mathrm{PT}_{\boldsymbol{\mathcal{A}}}]_{s}$ and $[Q]_{t}\in [\mathrm{PT}_{\boldsymbol{\mathcal{A}}}]_{t}$, and the following inequality holds
$$
\left([Q]_{t},t
\right)
<_{[
\mathbf{PT}_{\boldsymbol{\mathcal{A}}}
]}
\left([P]_{s},s
\right),
$$
where $\leq_{[
\mathbf{PT}_{\boldsymbol{\mathcal{A}}}
]}$ is the Artinian partial order on $\coprod[
\mathrm{PT}_{\boldsymbol{\mathcal{A}}}
]$ introduced in Definition~\ref{DPTQOrd}.

Note that, from Definition~\ref{DDCH}, the following chain of equalities holds
\begin{align*}
\mathrm{CH}^{(2)}_{t}\left(
\mathfrak{Q}^{(2)}
\right)&=
\eta^{(2,1)\sharp}_{t}\left(
\mathrm{CH}^{(1)\mathrm{m}}_{t}\left(
\mathrm{ip}^{([1],X)@}_{t}\left(
\left[
Q
\right]_{t}
\right)\right)\right)
\tag{1}
\\&=
\eta^{(2,1)\sharp}_{t}\left(
\mathrm{CH}^{(1)\mathrm{m}}_{t}\left(
\left[
\mathrm{ip}^{(1,X)@}_{t}\left(
Q
\right)\right]_{t}
\right)\right)
\tag{2}
\\&=
\eta^{(2,1)\sharp}_{t}\left(
\mathrm{CH}^{(1)}_{t}\left(
\mathrm{ip}^{(1,X)@}_{t}\left(
Q
\right)\right)\right).
\tag{3}
\end{align*}

The first equality unravels the definition of the second-order Curry-Howard mapping at a $(2,[1])$-identity second-order path; the second equality unravels the definition of the mapping $\mathrm{ip}^{([1],X)@}$ introduced in Definition~\ref{DPTQIp}; finally, the last equality unravels the definition of the monomorphic Curry-Howard mapping.

By a similar argument we also have that 
$$
\mathrm{CH}_{s}^{(2)}\left(
\mathfrak{P}^{(2)}
\right)
=
\eta^{(2,1)\sharp}_{s}\left(
\mathrm{CH}^{(1)}_{s}\left(
\mathrm{ip}^{(1,X)@}_{s}\left(
P
\right)\right)\right).
$$

Note, on the other hand, that according to Corollary~\ref{CPTQOrd} the inequality 
$$
\left(
\left[Q
\right]_{t},t
\right)
<_{[
\mathbf{PT}_{\boldsymbol{\mathcal{A}}}
]}
\left(\left[
P
\right]_{s},s
\right)$$
holds if, and only if, 
$
\mathrm{CH}^{(1)}_{t}(
\mathrm{ip}^{(1,X)@}_{t}(
Q
))
\leq_{\mathbf{T}_{\Sigma^{\boldsymbol{\mathcal{A}}}}(X)}
\mathrm{CH}^{(1)}_{s}(
\mathrm{ip}^{(1,X)@}_{s}(
P
)).
$

Since $\coprod\eta^{(2,[1])\sharp}$ is an order embedding of  $(\coprod\mathrm{T}_{\Sigma^{\boldsymbol{\mathcal{A}}}}(X), \leq_{\mathbf{T}_{\Sigma^{\boldsymbol{\mathcal{A}}}}(X)})$ into the partially ordered set $(\coprod\mathrm{T}_{\Sigma^{\boldsymbol{\mathcal{A}}^{(2)}}}(X), \leq_{\mathbf{T}_{\Sigma^{\boldsymbol{\mathcal{A}}^{(2)}}}(X)})$, we have that 
$$
\left(
\mathrm{CH}^{(2)}_{t}\left(
\mathfrak{Q}^{(2)}
\right),t
\right)
\leq_{\mathbf{T}_{\Sigma^{\boldsymbol{\mathcal{A}}^{(2)}}}(X)}
\left(
\mathrm{CH}^{(2)}_{s}
\left(
\mathfrak{P}^{(2)}
\right),s
\right).
$$

Case~(1) follows.

If~(2), i.e., if $\mathfrak{P}^{(2)}$ is a path of length strictly greater than one containing its first   second-order echelon at position $i\in\bb{\mathfrak{P}^{(2)}}$ then, depending on this position, we have the following subcases.

If $i=0$, then, according to Definition~\ref{DDOrd}, $\mathfrak{Q}^{(2)}$ is equal to $\mathfrak{P}^{(2),0,0}$ or $\mathfrak{P}^{(2),1,\bb{\mathfrak{P}^{(2)}}-1}$ and, according to Definition~\ref{DDCH}, we have that
$$
\mathrm{CH}^{(2)}_{s}
\left(
\mathfrak{P}^{(2)}
\right)
=
\mathrm{CH}_{s}^{(2)}\left(
\mathfrak{P}^{(2),1,\bb{\mathfrak{P}^{(2)}}-1}
\right)
\circ_{s}^{1\mathbf{T}_{\Sigma^{\boldsymbol{\mathcal{A}}}}(X)}
\mathrm{CH}_{s}^{(2)}
\left(
\mathfrak{P}^{(2),0,0}
\right).
$$
In any case, we have that 
$
(\mathrm{CH}^{(2)}_{s}(\mathfrak{Q}),s)
\leq_{\mathbf{T}_{\Sigma^{\boldsymbol{\mathcal{A}}^{(2)}}}(X)}
(\mathrm{CH}^{(2)}_{s}(\mathfrak{P}),s).
$

If $i>0$, then, according to Definition~\ref{DDOrd} $\mathfrak{Q}^{(2)}$ is equal to $\mathfrak{P}^{(2),0,i-1}$ or $\mathfrak{P}^{(2),i,\bb{\mathfrak{P}^{(2)}}-1}$ and, according to Definition~\ref{DDCH}, we have that
$$
\mathrm{CH}^{(2)}_{s}
\left(
\mathfrak{P}^{(2)}
\right)
=
\mathrm{CH}_{s}^{(2)}\left(
\mathfrak{P}^{(2),i,\bb{\mathfrak{P}^{(2)}}-1}
\right)
\circ_{s}^{1\mathbf{T}_{\Sigma^{\boldsymbol{\mathcal{A}}}}(X)}
\mathrm{CH}_{s}^{(2)}\left(
\mathfrak{P}^{(2),0,i-1}
\right).
$$
In any case, we have that 
$
(\mathrm{CH}^{(2)}_{s}(\mathfrak{Q}),s)
\leq_{\mathbf{T}_{\Sigma^{\boldsymbol{\mathcal{A}}^{(2)}}}(X)}
(\mathrm{CH}^{(2)}_{s}(\mathfrak{P}),s).
$

Case~(2) follows.

If~(3), i.e., if $\mathfrak{P}^{(2)}$ is an echelonless second-order path, then according to the different possibilities for $\mathfrak{P}^{(2)}$ we have the following subcases.

If $\mathfrak{P}^{(2)}$ is not head-constant, then, according to Definition~\ref{DDOrd}, for the maximum index $i\in\bb{\mathfrak{P}^{(2)}}$ for which $\mathfrak{P}^{(2),0,i}$ is a head-constant second-order path, then $\mathfrak{Q}^{(2)}$ is equal to $\mathfrak{P}^{(2),0,i}$ or $\mathfrak{P}^{(2), i+1,\bb{\mathfrak{P}^{(2)}}-1}$ and, according to Definition~\ref{DDCH}, we have that
$$
\mathrm{CH}^{(2)}_{s}\left(
\mathfrak{P}^{(2)}
\right)
=
\mathrm{CH}_{s}^{(2)}\left(
\mathfrak{P}^{(2),i+1,\bb{\mathfrak{P}^{(2)}}-1}
\right)
\circ_{s}^{1\mathbf{T}_{\Sigma^{\boldsymbol{\mathcal{A}}}}(X)}
\mathrm{CH}_{s}^{(2)}\left(
\mathfrak{P}^{(2),0,i}
\right).
$$
In any case, we conclude that 
$
(\mathrm{CH}^{(2)}_{s}(\mathfrak{Q}),s)
\leq_{\mathbf{T}_{\Sigma^{\boldsymbol{\mathcal{A}}^{(2)}}}(X)}
(\mathrm{CH}^{(2)}_{s}(\mathfrak{P}),s).
$

If $\mathfrak{P}^{(2)}$ is head-constant but not coherent then, according to Definition~\ref{DDOrd}, for the maximum index $i\in\bb{\mathfrak{P}^{(2)}}$ for which $\mathfrak{P}^{(2),0,i}$ is a coherent second-order path, then $\mathfrak{Q}^{(2)}$ is equal to $\mathfrak{P}^{(2),0,i}$ or $\mathfrak{P}^{(2), i+1,\bb{\mathfrak{P}^{(2)}}-1}$ and, according to Definition~\ref{DDCH}, we have that
$$
\mathrm{CH}^{(2)}_{s}
\left(
\mathfrak{P}^{(2)}
\right)
=
\mathrm{CH}_{s}^{(2)}
\left(
\mathfrak{P}^{(2),i+1,\bb{\mathfrak{P}^{(2)}}-1}
\right)
\circ_{s}^{1\mathbf{T}_{\Sigma^{\boldsymbol{\mathcal{A}}}}(X)}
\mathrm{CH}_{s}^{(2)}\left(
\mathfrak{P}^{(2),0,i}
\right).
$$
In any case, we conclude that 
$
(\mathrm{CH}^{(2)}_{s}(\mathfrak{Q}),s)
\leq_{\mathbf{T}_{\Sigma^{\boldsymbol{\mathcal{A}}^{(2)}}}(X)}
(\mathrm{CH}^{(2)}_{s}(\mathfrak{P}),s).
$

Therefore we are left with the case of $\mathfrak{P}^{(2)}$ being a head-constant coherent echelonless second-order path. Then according to Definition~\ref{DDHeadCt} there exists a unique word $\mathbf{s}\in S^{\star}-\{\lambda\}$, and a unique operation symbol $\tau\in \Sigma^{\boldsymbol{\mathcal{A}}}_{\mathbf{s},s}$ associated to $\mathfrak{P}^{(2)}$. Let $(\mathfrak{P}^{}_{j})_{j\in\bb{\mathbf{s}}}$ be the family of second-order paths in $\mathrm{Pth}_{\boldsymbol{\mathcal{A}}^{(2)}}$ we can extract from $\mathfrak{P}^{(2)}$ in virtue of Lemma~\ref{LDPthExtract}. Then, according to Definition~\ref{DDOrd}, $\mathfrak{Q}^{(2)}$ is equal to $\mathfrak{P}^{}_{j}$, for some $j\in\bb{\mathbf{s}}$, and, according to Definition~\ref{DDCH}, we have that 
$$
\mathrm{CH}^{(2)}_{s}
\left(
\mathfrak{P}^{(2)}
\right)=
\tau^{\mathbf{T}_{\Sigma^{\boldsymbol{\mathcal{A}}^{(2)}}}(X)}
\left(\left(
\mathrm{CH}^{(2)}_{s_{j}}\left(
\mathfrak{P}^{}_{j}
\right)\right)_{j\in\bb{\mathbf{s}}}
\right).
$$

In any case, we have that 
$
(\mathrm{CH}^{(2)}_{s_{j}}(\mathfrak{Q}^{(2)}),s_{j})
\leq_{\mathbf{T}_{\Sigma^{\boldsymbol{\mathcal{A}}^{(2)}}}(X)}
(\mathrm{CH}^{(2)}_{s}(\mathfrak{P}^{(2)}),s).
$

This concludes Case~(3).

This completes the base step.

\textsf{Inductive step of the induction.}

Assume the statement holds for sequences of length up to $m$, i.e., for every pair of sorts $t,s\in S$, if $(\mathfrak{Q}^{(2)},t), (\mathfrak{P}^{(2)},s)$ are pairs in $\coprod
\mathrm{Pth}_{\boldsymbol{\mathcal{A}}^{(2)}}$ such that  there exists a word $\mathbf{w}\in S^{\star}$ of length $\bb{\mathbf{w}}=m+1$ and a family of second-order paths $(\mathfrak{R}^{(2)}_{k})_{k\in\bb{\mathbf{w}}}$ in $\mathrm{Pth}_{\boldsymbol{\mathcal{A}}^{(2)},\mathbf{w}}$  such that $w_{0}=t$, $\mathfrak{R}^{(2)}_{0}=\mathfrak{Q}^{(2)}$, $w_{m}=s$, $\mathfrak{R}^{(2)}_{m}=\mathfrak{P}^{(2)}$ and 
for every $k\in m$, $(\mathfrak{R}^{(2)}_{k}, w_{k})\prec_{\mathbf{Pth}_{\boldsymbol{\mathcal{A}}^{(2)}}} (\mathfrak{R}^{(2)}_{k+1}, w_{k+1})$ then 
$
\mathrm{CH}^{(2)}_{t}(\mathfrak{Q}^{(2)})\in\mathrm{Subt}_{\Sigma^{\boldsymbol{\mathcal{A}}^{(2)}}}(\mathrm{CH}^{(2)}_{s}(\mathfrak{P}^{(2)}))_{t}.
$

Now we prove it for sequences of length $m+1$. 

Let $t,s$ be sorts in $S$, if $(\mathfrak{Q}^{(2)},t), (\mathfrak{P}^{(2)},s)$ are pairs in $\coprod
\mathrm{Pth}_{\boldsymbol{\mathcal{A}}^{(2)}}$ such that  there exists a word $\mathbf{w}\in S^{\star}$ of length $\bb{\mathbf{w}}=m+2$ and a family of second-order paths $(\mathfrak{R}^{(2)}_{k})_{k\in\bb{\mathbf{w}}}$ in $\mathrm{Pth}_{\boldsymbol{\mathcal{A}}^{(2)},\mathbf{w}}$  such that $w_{0}=t$, $\mathfrak{R}^{(2)}_{0}=\mathfrak{Q}^{(2)}$, $w_{m+1}=s$, $\mathfrak{R}^{(2)}_{m+1}=\mathfrak{P}^{(2)}$ and 
for every $k\in m+1$, $(\mathfrak{R}^{(2)}_{k}, w_{k})\prec_{\mathbf{Pth}_{\boldsymbol{\mathcal{A}}^{(2)}}} (\mathfrak{R}^{(2)}_{k+1}, w_{k+1})$.

Consider the word $\mathbf{w}^{0,m}$ of length $\bb{\mathbf{w}^{0,m}}=m+1$ and the family of second-order paths $(\mathfrak{R}^{(2)}_{k})_{k\in\bb{\mathbf{w}^{0,m}}}$ in $\mathrm{Pth}_{\boldsymbol{\mathcal{A}}^{(2)},\mathbf{w}^{0,m}}$. This is a sequence of length $m$ instantiating that 
$
(\mathfrak{Q}^{(2)},t)
\leq_{\mathbf{Pth}_{\boldsymbol{\mathcal{A}}^{(2)}}}
(\mathfrak{R}^{(2)}_{m}, w_{m})$. By the inductive hypothesis, we have that 
$\mathrm{CH}^{(2)}_{t}(\mathfrak{Q}^{(2)})
\in\mathrm{Subt}_{\Sigma^{\boldsymbol{\mathcal{A}}^{(2)}}}(\mathrm{CH}^{(2)}_{w_{m}}(\mathfrak{R}^{(2)}_{m}))_{t}.
$

Now, consider the sequence of second-order paths $(\mathfrak{R}^{(2)}_{m},\mathfrak{P}^{(2)})$ it is a one-step sequence of second-order paths in $\mathrm{Pth}_{\boldsymbol{\mathcal{A}}^{(2)},\mathbf{w}^{m,m+1}}$ satisfying that $(\mathfrak{R}^{(2)}_{m},w_{m})\prec_{\mathbf{Pth}_{\boldsymbol{\mathcal{A}}^{(2)}}} (\mathfrak{P}^{(2)}, s)$. By the base case, we have that 
$
\mathrm{CH}^{(2)}_{w_{m}}(\mathfrak{R}^{(2)}_{m})
\in\mathrm{Subt}_{\Sigma^{\boldsymbol{\mathcal{A}}^{(2)}}}(\mathrm{CH}^{(2)}_{s}(\mathfrak{P}^{(2)}))_{w_{m}}.
$

Hence, by the transitivity of the partial order $\leq_{\mathbf{T}_{\Sigma^{\boldsymbol{\mathcal{A}}^{(2)}}}(X)}$, we have that 
$$\mathrm{CH}^{(2)}_{t}
\left(\mathfrak{Q}^{(2)}
\right)
\in\mathrm{Subt}_{\Sigma^{\boldsymbol{\mathcal{A}}^{(2)}}}
\left(\mathrm{CH}^{(2)}_{s}
\left(\mathfrak{P}^{(2)}
\right)\right)_{t}.
$$

This finishes the proof.
\end{proof}

\chapter{
\texorpdfstring
{On the kernel of the second-order Curry-Howard mapping}
{The kernel of the second-order Curry-Howard mapping}
}\label{S2E}

This chapter is devoted to the study of the second-order Curry-Howard mapping. We first prove that the second-order Curry-Howard mapping is able to completely characterise $(2,[1])$-identity second-order paths. In fact, we prove that $(2,[1])$-identity second-order paths in the kernel of $\mathrm{CH}^{(2)}$ must be necessarily equal. In an analogous  manner, the second-order Curry-Howard mapping is able to completely characterise $(2,0)$-identity second-order paths. In fact, we prove that $(2,0)$-identity second-order paths in the kernel of $\mathrm{CH}^{(2)}$ must be necessarily equal. The second-order Curry-Howard mapping is also capable of characterising second-order echelons. We prove that second-order paths in the kernel of the second-order Curry-Howard mapping necessarily have the same length, the same $([1],2)$-source, the same $([1],2)$-target, the same $(0,2)$-source and the same $(0,2)$-target. This, in fact, helps to demonstrate that second-order paths of length one are completely characterised by the second-order Curry-Howard mapping. Indeed, two second-order paths of length one in  the kernel of $\mathrm{CH}^{(2)}$ must be necessarily equal. Moreover, the second-order paths in the kernel of the second-order Curry-Howard mapping preserve the coherence.  The most important result of this chapter is the one that proves that $\mathrm{Ker}(\mathrm{CH}^{(2)})$ is a closed $\Sigma^{\boldsymbol{\mathcal{A}}^{(2)}}$-congruence on $\mathbf{Pth}_{\boldsymbol{\mathcal{A}}^{(2)}}$.


First of all, we identify $\mathcal{A}^{(2)}$ with a subset of terms in $\mathrm{T}_{\Sigma^{\boldsymbol{\mathcal{A}}^{(2)}}}(X)$ by means of its image under the $S$-sorted mapping $\eta^{(2,\mathcal{A}^{(2)})}$ presented in Definition~\ref{DDEta}, that is 
$$
\eta^{(2,\mathcal{A}^{(2)})}\left[
\mathcal{A}^{(2)}
\right]=
\left(
\left\lbrace
\mathfrak{p}^{(2)\mathbf{T}_{\Sigma^{\boldsymbol{\mathcal{A}}^{(2)}}}(X)}\in 
\mathrm{T}_{\Sigma^{\boldsymbol{\mathcal{A}}^{(2)}}}(X)_{s}
\,\middle|\,
\mathfrak{p}^{(2)}\in\mathcal{A}^{(2)}_{s}
\right\rbrace
\right)
_{s\in S}.
$$

We begin by showing that the value of the second-order Curry-Howard mapping at a non-$(2,[1])$-identity second-order path always has a subterm which contains a second-order rewrite rule (of some sort) in $\mathcal{A}^{(2)}$. 

\begin{proposition}\label{PDCHRewId}
Let $s$ be a sort in $S$ and $\mathfrak{P}^{(2)}$ a non-$(2,[1])$-identity second-order path in $\mathrm{Pth}_{\boldsymbol{\mathcal{A}}^{(2)},s}$. Then 
$$\eta^{(2,\mathcal{A}^{(2)})}\left[
\mathcal{A}^{(2)}
\right]\cap\mathrm{Subt}_{\Sigma^{\boldsymbol{\mathcal{A}}^{(2)}}}
\left(
\mathrm{CH}^{(2)}_{s}\left(
\mathfrak{P}^{(2)}
\right)\right)
\neq\varnothing^{S}.$$
\end{proposition}
\begin{proof}
We prove the statement by Artinian induction on $(\coprod\mathrm{Pth}_{\boldsymbol{\mathcal{A}}^{(2)}}, \leq_{\mathbf{Pth}_{\boldsymbol{\mathcal{A}}^{(2)}}})$.

\textsf{Base step of the Artinian induction}

Let $(\mathfrak{P}^{(2)},s)$ be a minimal element of $(\coprod\mathrm{Pth}_{\boldsymbol{\mathcal{A}}^{(2)}}, \leq_{\mathbf{Pth}_{\boldsymbol{\mathcal{A}}^{(2)}}})$. Thus, by Proposition~\ref{PDMinimal}, the second-order path $\mathfrak{P}^{(2)}$ is either (1) a $(2,[1])$-identity second-order path on a minimal path term class or (2) a   second-order echelon. Since we are assuming that $\mathfrak{P}^{(2)}$ is a non-$(2,[1])$-identity second-order path, we have that $\mathfrak{P}^{(2)}$ is a   second-order echelon. Hence 
$$
\mathrm{CH}^{(2)}_{s}
\left(
\mathfrak{P}^{(2)}
\right)=\mathfrak{p}^{(2)\mathbf{T}_{\Sigma^{\boldsymbol{\mathcal{A}}^{(2)}}}(X)},
$$
for $\mathfrak{p}^{(2)}\in\mathcal{A}^{(2)}_{s}$, the unique second-order rewrite rule occurring in $\mathfrak{P}^{(2)}$.

Therefore
$$
\mathfrak{p}^{(2)\mathbf{T}_{\Sigma^{\boldsymbol{\mathcal{A}}^{(2)}}}(X)}
\in 
\eta^{(2,\mathcal{A}^{(2)})}
\left[
\mathcal{A}^{(2)}
\right]_{s}
\cap
\mathrm{Subt}_{\Sigma^{\boldsymbol{\mathcal{A}}^{(2)}}}\left(
\mathrm{CH}^{(2)}_{s}\left(
\mathfrak{P}^{(2)}
\right)\right)_{s}.
$$

\textsf{Inductive step of the Artinian induction}

Let $(\mathfrak{P}^{(2)},s)$ be a non-minimal element of $(\coprod\mathrm{Pth}_{\boldsymbol{\mathcal{A}}^{(2)}}, \leq_{\mathbf{Pth}_{\boldsymbol{\mathcal{A}}^{(2)}}})$ satisfying that $\mathfrak{P}^{(2)}$ is a non-$(2,[1])$-identity second-order path. Let us suppose that, for every sort $t\in S$ and every second-order path $\mathfrak{Q}^{(2)}$ in $\mathrm{Pth}_{\boldsymbol{\mathcal{A}}^{(2)},t}$, if $(\mathfrak{Q}^{(2)},t)$ $\leq_{\mathbf{Pth}_{\boldsymbol{\mathcal{A}}^{(2)}}}$-precedes $(\mathfrak{P}^{(2)},s)$, then the statement holds for $\mathfrak{Q}^{(2)}$, i.e., if $\mathfrak{Q}^{(2)}$ is a non-$(2,[1])$-identity second-order path in  $\mathrm{Pth}_{\boldsymbol{\mathcal{A}}^{(2)},t}$, then
$$
\eta^{(2,\mathcal{A}^{(2)})}\left[
\mathcal{A}^{(2)}
\right]
\cap
\mathrm{Subt}_{\Sigma^{\boldsymbol{\mathcal{A}}^{(2)}}}\left(
\mathrm{CH}^{(2)}_{t}\left(
\mathfrak{Q}^{(2)}
\right)\right)
\neq\varnothing^{S}.
$$

Since $(\mathfrak{P}^{(2)},s)$ is a non-minimal element of $(\coprod\mathrm{Pth}_{\boldsymbol{\mathcal{A}}^{(2)}}, \leq_{\mathbf{Pth}_{\boldsymbol{\mathcal{A}}^{(2)}}})$ and taking into account that $\mathfrak{P}^{(2)}$ is a non-$(2,[1])$-identity second-order path, we have, by Lemma~\ref{LDOrdI}, that $\mathfrak{P}^{(2)}$ is either (1) a second-order path of length strictly greater than one containing at least one   second-order echelon, or (2) an echelonless second-order path.

If~(1), then let $i\in\bb{\mathfrak{P}^{(2)}}$ be the first index for which the one-step subpath $\mathfrak{P}^{(2),i,i}$ is a   second-order echelon. In this case, according to the nature of $i$, we have that, if $i=0$, then $\mathfrak{P}^{(2),0,0}$ is a   second-order echelon and, if $i>0$, then $\mathfrak{P}^{(2),i,\bb{\mathfrak{P}^{(2)}}-1}$ contains a   second-order echelon.

Note that in any case we have found a subpath $\mathfrak{Q}^{(2)}$ of $\mathfrak{P}^{(2)}$ which is a non-$(2,[1])$-identity second-order path and such that $(\mathfrak{Q}^{(2)},s)$ $\leq_{\mathbf{Pth}_{\boldsymbol{\mathcal{A}}^{(2)}}}$-precedes $(\mathfrak{P}^{(2)},s)$. Then, by induction, we have that 
$
\eta^{(2,\mathcal{A}^{(2)})}[\mathcal{A}^{(2)}]
\cap
\mathrm{Subt}_{\Sigma^{\boldsymbol{\mathcal{A}}^{(2)}}}(
\mathrm{CH}^{(2)}_{s}(
\mathfrak{Q}^{(2)}
))\neq\varnothing^{S}.
$
By Proposition~\ref{PDCHMono}, we also have that $\mathrm{CH}^{(2)}_{s}(\mathfrak{Q}^{(2)})\in\mathrm{Subt}_{\Sigma^{\boldsymbol{\mathcal{A}}^{(2)}}}(\mathrm{CH}^{(2)}_{s}(\mathfrak{P}^{(2)}))_{s}$, thus we have that 
$
\eta^{(2,\mathcal{A}^{(2)})}[\mathcal{A}^{(2)}]
\cap
\mathrm{Subt}_{\Sigma^{\boldsymbol{\mathcal{A}}^{(2)}}}(
\mathrm{CH}^{(2)}_{s}(
\mathfrak{P}^{(2)}
))\neq\varnothing^{S}.
$

If~(2), i.e., if $\mathfrak{P}^{(2)}$ is an echelonless second-order path, it could be the case that (2.1) $\mathfrak{P}^{(2)}$ is an echelonless second-order path that is not head-constant, or (2.2) $\mathfrak{P}^{(2)}$ is a head-constant echelonless second-order path that is not coherent or (2.3) $\mathfrak{P}^{(2)}$ is a head-constant coherent echelonless second-order path.

If~(2.1), let $i\in\bb{\mathfrak{P}^{(2)}}$ be the greatest index $i\in\bb{\mathfrak{P}^{(2)}}$ for which $\mathfrak{P}^{(2),0,i}$ is a head-constant second-order path. Note that $\mathfrak{P}^{(2),0,i}$ is a non-$(2,[1])$-identity second-order path such that $(\mathfrak{P}^{(2),0,i},s)$  $\leq_{\mathbf{Pth}_{\boldsymbol{\mathcal{A}}^{(2)}}}$-precedes $(\mathfrak{P}^{(2)},s)$. Then, by induction, we have that 
$
\eta^{(2,\mathcal{A}^{(2)})}[\mathcal{A}^{(2)}]
\cap
\mathrm{Subt}_{\Sigma^{\boldsymbol{\mathcal{A}}^{(2)}}}(
\mathrm{CH}^{(2)}_{s}(
\mathfrak{P}^{(2),0,i}
))\neq\varnothing^{S}.
$
By Proposition~\ref{PDCHMono}, we also have that $\mathrm{CH}^{(2)}_{s}(\mathfrak{P}^{(2),0,i})\in\mathrm{Subt}_{\Sigma^{\boldsymbol{\mathcal{A}}^{(2)}}}(\mathrm{CH}^{(2)}_{s}(\mathfrak{P}^{(2)}))_{s}$, thus we have that 
$
\eta^{(2,\mathcal{A}^{(2)})}[\mathcal{A}^{(2)}]
\cap
\mathrm{Subt}_{\Sigma^{\boldsymbol{\mathcal{A}}^{(2)}}}(
\mathrm{CH}^{(2)}_{s}(
\mathfrak{P}^{(2)}
))\neq\varnothing^{S}.
$

If~(2.2), let $i\in\bb{\mathfrak{P}^{(2)}}$ be the greatest index $i\in\bb{\mathfrak{P}^{(2)}}$ for which $\mathfrak{P}^{(2),0,i}$ is a coherent second-order path. Note that $\mathfrak{P}^{(2),0,i}$ is a non-$(2,[1])$-identity second-order path such that $(\mathfrak{P}^{(2),0,i},s)$  $\leq_{\mathbf{Pth}_{\boldsymbol{\mathcal{A}}^{(2)}}}$-precedes $(\mathfrak{P}^{(2)},s)$. Then, by induction, we have that 
$
\eta^{(2,\mathcal{A}^{(2)})}[\mathcal{A}^{(2)}]
\cap
\mathrm{Subt}_{\Sigma^{\boldsymbol{\mathcal{A}}^{(2)}}}(
\mathrm{CH}^{(2)}_{s}(
\mathfrak{P}^{(2),0,i}
))\neq\varnothing^{S}.
$
By Proposition~\ref{PDCHMono}, we also have that $\mathrm{CH}^{(2)}_{s}(\mathfrak{P}^{(2),0,i})\in\mathrm{Subt}_{\Sigma^{\boldsymbol{\mathcal{A}}^{(2)}}}(\mathrm{CH}^{(2)}_{s}(\mathfrak{P}^{(2)}))_{s}$, thus we have that 
$
\eta^{(2,\mathcal{A}^{(2)})}[\mathcal{A}^{(2)}]
\cap
\mathrm{Subt}_{\Sigma^{\boldsymbol{\mathcal{A}}^{(2)}}}(
\mathrm{CH}^{(2)}_{s}(
\mathfrak{P}^{(2)}
))\neq\varnothing^{S}.
$

If~(2.3), then in virtue of Definition~\ref{DDHeadCt} there exists a unique word $\mathbf{s}\in S^{\star}-\{\lambda\}$ and a unique operation symbol $\tau\in\Sigma^{\boldsymbol{\mathcal{A}}}_{\mathbf{s},s}$ associated to $\mathfrak{P}$. Let $(\mathfrak{P}^{(2)}_{j})_{j\in\bb{\mathbf{s}}}$ be the family of second-order paths in $\mathrm{Pth}_{\boldsymbol{\mathcal{A}}^{(2)},\mathbf{s}}$ which, in virtue of Lemma~\ref{LDPthExtract}, we can extract from $\mathfrak{P}^{(2)}$. Since $\bb{\mathfrak{P}^{(1)}}\geq 1$, there exists at least one $j\in\bb{\mathbf{s}}$ such that $\mathfrak{P}^{(2)}_{j}$ is a non-$(2,[1])$-identity second-order path and, in addition, $(\mathfrak{P}^{(2)}_{j},s_{j})$ $\leq_{\mathbf{Pth}_{\boldsymbol{\mathcal{A}}^{(2)}}}$-precedes $(\mathfrak{P}^{(2)},s)$. Then, by induction, we have that 
$
\eta^{(2,\mathcal{A}^{(2)})}[\mathcal{A}^{(2)}]
\cap
\mathrm{Subt}_{\Sigma^{\boldsymbol{\mathcal{A}}^{(2)}}}(
\mathrm{CH}^{(2)}_{s_{j}}(
\mathfrak{P}^{(2)}_{j}
))\neq\varnothing^{S}.
$
By Proposition~\ref{PDCHMono}, we also have that $\mathrm{CH}^{(2)}_{s_{j}}(\mathfrak{P}^{(2)}_{j})\in\mathrm{Subt}_{\Sigma^{\boldsymbol{\mathcal{A}}^{(2)}}}(\mathrm{CH}^{(2)}_{s}(\mathfrak{P}^{(2)}))_{s_{j}}$, thus we have that 
$
\eta^{(2,\mathcal{A}^{(2)})}[\mathcal{A}^{(2)}]
\cap
\mathrm{Subt}_{\Sigma^{\boldsymbol{\mathcal{A}}^{(2)}}}(
\mathrm{CH}^{(2)}_{s}(
\mathfrak{P}^{(2)}
))\neq\varnothing^{S}.
$

This finishes the proof.
\end{proof}

In the following corollary we characterize the $(2,[1])$-identity second-order paths with the help of the second-order Curry-Howard mapping.

\begin{restatable}{corollary}{CDCHId}
\label{CDCHId} Let $s$ be a sort in $S$ and $\mathfrak{P}^{(2)}$ a second-order path in $\mathrm{Pth}_{\boldsymbol{\mathcal{A}}^{(2)},s}$. The following statements are equivalent.
\begin{enumerate}
\item[(i)] $\mathfrak{P}^{(2)}$ is a $(2,[1])$-identity second-order path.
\item[(ii)] $\mathrm{CH}^{(2)}_{s}(\mathfrak{P}^{(2)})=\eta^{(2,1)\sharp}_{s}(\mathrm{CH}^{(1)\mathrm{m}}_{s}(\mathrm{ip}^{([1],X)@}_{s}(\mathrm{sc}^{([1],2)}_{s}(\mathfrak{P}))))$.
\item[(iii)] $\mathrm{CH}^{(2)}_{s}(\mathfrak{P}^{(2)})=\eta^{(2,1)\sharp}_{s}(\mathrm{CH}^{(1)\mathrm{m}}_{s}(\mathrm{ip}^{([1],X)@}_{s}(\mathrm{tg}^{([1],2)}_{s}(\mathfrak{P}))))$.
\item[(iv)] $\mathrm{CH}^{(2)}_{s}(\mathfrak{P}^{(2)})\in \eta^{(2,1)\sharp}[\mathrm{PT}_{\boldsymbol{\mathcal{A}}}]_{s}$.
\end{enumerate}
\end{restatable}
\begin{proof}
We first prove that (i) implies (ii). Let $\mathfrak{P}^{(2)}$ be a $(2,[1])$-identity second-order path of sort $s\in S$ then, for some path term class $[P]_{s}$ in $[\mathrm{PT}_{\boldsymbol{\mathcal{A}}}]_{s}$, we have that $\mathfrak{P}^{(2)}=\mathrm{ip}^{(2,[1])\sharp}_{s}([P]_{s})$.

The following chain of equalities holds
\allowdisplaybreaks
\begin{align*}
\mathrm{CH}^{(2)}_{s}\left(
\mathfrak{P}^{(2)}
\right)&=
\mathrm{CH}^{(2)}_{s}\left(
\mathrm{ip}^{(2,[1])\sharp}_{s}\left(
\left[P
\right]_{s}
\right)\right)
\tag{1}
\\&=
\eta^{(2,1)\sharp}_{s}\left(
\mathrm{CH}^{(1)\mathrm{m}}_{s}\left(
\mathrm{ip}^{([1],X)@}_{s}\left(
\left[P
\right]_{s}
\right)\right)\right)
\tag{2}
\\&=
\eta^{(2,1)\sharp}_{s}\left(
\mathrm{CH}^{(1)\mathrm{m}}_{s}\left(
\mathrm{ip}^{([1],X)@}_{s}\left(
\mathrm{sc}^{([1],2)}_{s}\left(
\mathfrak{P}^{(2)}
\right)\right)\right)\right).
\tag{3}
\end{align*}

The first equality unravels the description of $\mathfrak{P}^{(2)}$ as a $(2,[1])$-identity second-order path; the second equality follows from Proposition~\ref{PDCHDUId}; finally, the last equality recovers the $([1],2)$-source of the second-order path $\mathfrak{P}^{(2)}$.

That (ii) implies (iv) is straightforward, since, by Lemmas~\ref{LWCong} and~\ref{LThetaCong}, $\mathrm{CH}^{(1)\mathrm{m}}_{s}(\mathrm{ip}^{([1],X)@}_{s}(\mathrm{sc}^{([1],2)}_{s}(\mathfrak{P}^{(2)})))$ is a term in $\mathrm{PT}_{\boldsymbol{\mathcal{A}}, s}$.

Finally, we prove that (iv) implies (i). By contraposition, assume that $\mathfrak{P}^{(2)}$ is a non-$(2,[1])$-identity second-order path. Then, by Proposition~\ref{PDCHRewId},
$$
\eta^{(2,\mathcal{A}^{(2)})}\left[
\mathcal{A}^{(2)}
\right]
\cap\mathrm{Subt}_{\Sigma^{\boldsymbol{\mathcal{A}}^{(2)}}}\left(
\mathrm{CH}^{(2)}_{s}\left(
\mathfrak{P}^{(2)}
\right)\right)
\neq
\varnothing^{S},
$$
contradicting the assumption that the term $\mathrm{CH}^{(2)}_{s}(\mathfrak{P}^{(2)})$ belongs to $\eta^{(2,1)\sharp}[\mathrm{PT}_{\boldsymbol{\mathcal{A}}}]_{s}$.

Replacing the $([1],2)$-source by the $([1],2)$-target in the previous development follows easily since they coincide on $(2,[1])$-identity second-order paths.

This finishes the proof.
\end{proof}

We can now complete the characterization of the non-$(2,[1])$-identity second-order paths.

\begin{corollary}\label{CDCHRew}
Let $s$ be a sort in $S$ and $\mathfrak{P}^{(2)}$ a second-order path in $\mathrm{Pth}_{\boldsymbol{\mathcal{A}}^{(2)},s}$. The following statements are equivalent.
\begin{enumerate}
\item[(i)] $\mathfrak{P}^{(2)}$ is a  non-$(2,[1])$-identity second-order path.
\item[(ii)] $\eta^{(2,\mathcal{A}^{(2)})}[\mathcal{A}^{(2)}]\cap\mathrm{Subt}_{\Sigma^{\boldsymbol{\mathcal{A}}^{(2)}}}(\mathrm{CH}^{(2)}_{s}(\mathfrak{P}^{(2)}))\neq\varnothing^{S}$.
\end{enumerate}
\end{corollary}

\begin{proof}
That (i) implies (ii) was proved in Proposition~\ref{PDCHRewId}.

Now we prove that (ii) implies (i). By contraposition, let us assume that $\mathfrak{P}^{(2)}$ is a $(2,[1])$-identity second-order path. Then, by Corollary~\ref{CDCHId}, we have that $\mathrm{CH}^{(2)}_{s}(\mathfrak{P})$ is a term in $\eta^{(2,1)\sharp}[\mathrm{PT}_{\boldsymbol{\mathcal{A}}}]_{s}$, thus it cannot contain subterms of any sort from $\eta^{(2,\mathcal{A}^{(2)})}[\mathcal{A}^{(2)}]$.
\end{proof}

Another interesting consequence of Corollary~\ref{CDCHId} is the following one. It states that a second-order path that has the same image under the second-order Curry-Howard mapping than that of a $(2,[1])$-identity second-order path implies that both second-order paths are  equal.

\begin{restatable}{corollary}{CDCHUId}
\label{CDCHUId}
Let $s$ be a sort in $S$ and $(\mathfrak{P}^{(2)},\mathfrak{Q}^{(2)})\in\mathrm{Ker}(\mathrm{CH}^{(2)})_{s}$. If $\mathfrak{P}^{(2)}$ or $\mathfrak{Q}^{(2)}$ is a $(2,[1])$-identity second-order path, then $\mathfrak{P}^{(2)}=\mathfrak{Q}^{(2)}$.
\end{restatable}
\begin{proof}
Let us assume that $\mathfrak{P}^{(2)}$ is a $(2,[1])$-identity second-order path. Then there exists a path term class $[P]_{s}$ in $[\mathrm{PT}_{\boldsymbol{\mathcal{A}}}]_{s}$ for which 
$\mathfrak{P}^{(2)}=\mathrm{ip}^{(2,[1])\sharp}_{s}([P]_{s})$. 

The following chain of equalities holds 
\allowdisplaybreaks
\begin{align*}
\mathrm{CH}^{(2)}_{s}
\left(
\mathfrak{P}^{(2)}
\right)
&=
\mathrm{CH}^{(2)}_{s}
\left(
\mathrm{ip}^{(2,[1])\sharp}_{s}\left(
[P]_{s}
\right)
\right)
\tag{1}
\\&=
\eta^{(2,1)\sharp}_{s}\left(
\mathrm{CH}^{(1)\mathrm{m}}_{s}\left(
\mathrm{ip}^{([1],X)@}_{s}\left(
[P]_{s}
\right)\right)\right)
\tag{2}
\\&=
\eta^{(2,1)\sharp}_{s}\left(
\mathrm{CH}^{(1)}_{s}\left(
\mathrm{ip}^{(1,X)@}_{s}\left(
P
\right)\right)\right)
.
\tag{3}
\end{align*}

The first equality unravels the description of $\mathfrak{P}^{(2)}$ as a $(2,[1])$-identity second-order path; the second equality follows from Corollary~\ref{CDCHId}; finally, the last equality follows by taking into account the definitions of the mappings $\mathrm{ip}^{(1,X)@}$ and $\mathrm{CH}^{(1)\mathrm{m}}$.

But, since $(\mathfrak{P}^{(2)},\mathfrak{Q}^{(2)})\in\mathrm{Ker}(\mathrm{CH}^{(2)})_{s}$, we have that $\mathrm{CH}^{(2)}_{s}(\mathfrak{Q}^{(2)})$ is a term in $\eta^{(2,1)\sharp}[\mathrm{PT}_{\boldsymbol{\mathcal{A}}}]_{s}$. Hence, by Corollary~\ref{CDCHId}, it follows that  $\mathfrak{Q}^{(2)}$ is a $(2,[1])$-identity second-order paths. Therefore, there exists some path term class $[Q]_{s}$ in $[\mathrm{PT}_{\boldsymbol{\mathcal{A}}}]_{s}$ for which 
$\mathfrak{Q}^{(2)}=\mathrm{ip}^{(2,[1])\sharp}_{s}([Q]_{s})$. 

By a similar computation as the one above, we have that  
$$\mathrm{CH}^{(2)}_{s}\left(
\mathfrak{Q}^{(2)}
\right)
=
\eta^{(2,1)\sharp}_{s}\left(
\mathrm{CH}^{(1)}_{s}\left(
\mathrm{ip}^{(1,X)@}_{s}\left(
Q
\right)\right)\right)
.$$ 

Taking into account that $(\mathfrak{P}^{(2)},\mathfrak{Q}^{(2)})\in\mathrm{Ker}(\mathrm{CH}^{(2)})_{s}$ and the fact that, by Proposition~\ref{PDUEmb} $\eta^{(2,1)\sharp}$ is injective, we have that 
$$
\mathrm{CH}^{(1)}_{s}\left(
\mathrm{ip}^{(1,X)@}_{s}\left(
P
\right)\right)=
\mathrm{CH}^{(1)}_{s}\left(
\mathrm{ip}^{(1,X)@}_{s}\left(
Q
\right)\right).
$$

Note that, by Lemma~\ref{LWCong}, the following pairs are related under $\Theta^{[1]}$
$$
\left(
P, \mathrm{CH}^{(1)}_{s}\left(
\mathrm{ip}^{(1,X)@}_{s}\left(
P
\right)\right) \right), 
\left(Q, \mathrm{CH}^{(1)}_{s}\left(
\mathrm{ip}^{(1,X)@}_{s}\left(
Q
\right)\right) \right)
\in\Theta^{[1]}_{s};
$$

By the transitivity of $\Theta^{[1]}$ we have that 
$
(P,Q)\in\Theta^{[1]}_{s}.
$
Thus, we can affirm that
$$
\mathfrak{P}^{(2)}=\mathrm{ip}^{(2,[1])\sharp}_{s}\left(
[P]_{s}
\right)
=
\mathrm{ip}^{(2,[1])\sharp}_{s}\left(
[Q]_{s}
\right)
=
\mathfrak{Q}^{(2)}.
$$

The statement follows.
\end{proof}

In the following corollary, we characterize the $(2,0)$-identity second-order paths in terms of their image under the second-order Curry-Howard mapping.

\begin{restatable}{corollary}{CDCHZId}
\label{CDCHZId}
Let $s$ be a sort in $S$ and $\mathfrak{P}^{(2)}$ a second-order path in $\mathrm{Pth}_{\boldsymbol{\mathcal{A}}^{(2)},s}$. The following statements are equivalent.
\begin{enumerate}
\item[(i)] $\mathfrak{P}^{(2)}$ is a $(2,0)$-identity second-order path.
\item[(ii)] $\mathrm{CH}^{(2)}_{s}(\mathfrak{P}^{(2)})=\eta^{(2,0)\sharp}_{s}(\mathrm{sc}^{(0,2)}_{s}(\mathfrak{P}^{(2)}))$.
\item[(iii)] $\mathrm{CH}^{(2)}_{s}(\mathfrak{P}^{(2)})=\eta^{(2,0)\sharp}_{s}(\mathrm{tg}^{(0,2)}_{s}(\mathfrak{P}^{(2)}))$.
\item[(iv)] $\mathrm{CH}^{(2)}_{s}(\mathfrak{P}^{(2)})\in\eta^{(2,0)\sharp}[\mathrm{T}_{\Sigma}(X)]_{s}$.
\end{enumerate}
\end{restatable}
\begin{proof}
We first prove that (i) implies (ii). Let $\mathfrak{P}^{(2)}$ be a $(2,0)$-identity second-order path of sort $s$. Then, for some term $P\in\mathrm{T}_{\Sigma}(X)_{s}$ we have that $\mathfrak{P}^{(2)}=\mathrm{ip}^{(2,0)\sharp}_{s}(P)$.

The following chain of equalities holds
\allowdisplaybreaks
\begin{align*}
\mathrm{CH}^{(2)}_{s}\left(
\mathfrak{P}^{(2)}
\right)
&=
\mathrm{CH}^{(2)}_{s}
\left(\mathrm{ip}^{(2,0)\sharp}_{s}(P)
\right)
\tag{1}
\\&=
\eta^{(2,0)\sharp}_{s}(P)
\tag{2}
\\&=
\eta^{(2,0)\sharp}_{s}\left(
\mathrm{sc}^{(0,2)}_{s}\left(
\mathfrak{P}^{(2)}
\right)\right).
\tag{3}
\end{align*}

The first equality unravels the description of $\mathfrak{P}^{(2)}$ as a $(2,0)$-identity second-order path; the second equality follows from Proposition~\ref{PDCHDZId}; finally, the last equality recovers the $(0,2)$-source of the second-order path $\mathfrak{P}^{(2)}$.

That (ii) implies (iv) is straightforward since $\mathrm{sc}^{(0,2)}_{s}(\mathfrak{P}^{(2)})$ is a term in $\mathrm{T}_{\Sigma}(X)_{s}$.

Assume (iv), i.e., that $\mathrm{CH}^{(2)}_{s}(\mathfrak{P}^{(2)})$ is a term in $\eta^{(2,0)\sharp}_{s}[\mathrm{T}_{\Sigma}(X)]_{s}$ then, by Proposition~\ref{PEmb}, since
$
\eta^{(1,0)\sharp}[\mathrm{T}_{\Sigma}(X)]
\subseteq
\mathrm{PT}_{\boldsymbol{\mathcal{A}}},
$
we have that
$
\eta^{(2,1)\sharp}[
\eta^{(1,0)\sharp}[
\mathrm{T}_{\Sigma}(X)
]]
\subseteq
\eta^{(2,1)\sharp}[
\mathrm{PT}_{\boldsymbol{\mathcal{A}}}
].
$
Hence, by Proposition~\ref{PDEmb}, we have that 
$
\eta^{(2,0)\sharp}[
\mathrm{T}_{\Sigma}(X)
]
\subseteq
\eta^{(2,1)\sharp}[
\mathrm{PT}_{\boldsymbol{\mathcal{A}}}].
$ 
Therefore, $\mathrm{CH}^{(2)}_{s}(\mathfrak{P}^{(2)})\in\eta^{(2,1)\sharp}[\mathrm{PT}_{\boldsymbol{\mathcal{A}}}]_{s}$.

By Corollary~\ref{CDCHId}, we have that $\mathfrak{P}^{(2)}$ is a $(2,[1])$-identity second-order path. Hence there exists a path term class $[P]_{s}$ in $[\mathrm{PT}_{\boldsymbol{\mathcal{A}}}]_{s}$ for which $\mathfrak{P}^{(2)}=\mathrm{ip}^{(2,[1])\sharp}_{s}([P]_{s})$.

The following chain of equalities holds
\allowdisplaybreaks
\begin{align*}
\mathrm{CH}^{(2)}_{s}\left(
\mathfrak{P}^{(2)}
\right)&=
\mathrm{CH}^{(2)}_{s}\left(
\mathrm{ip}^{(2,[1])\sharp}_{s}\left(
[P
]_{s}
\right)\right)
\tag{1}
\\&=
\eta^{(2,1)\sharp}_{s}\left(
\mathrm{CH}^{(1)\mathrm{m}}_{s}\left(
\mathrm{ip}^{([1],X)@}_{s}\left(
[P
]_{s}
\right)\right)\right)
\tag{2}
\\&=
\eta^{(2,1)\sharp}_{s}\left(
\mathrm{CH}^{(1)\mathrm{m}}_{s}\left(
\left[\mathrm{ip}^{(1,X)@}_{s}\left(
P
\right)\right]_{s}
\right)\right)
\tag{3}
\\&=
\eta^{(2,1)\sharp}_{s}\left(
\mathrm{CH}^{(1)}_{s}\left(
\mathrm{ip}^{(1,X)@}_{s}\left(
P
\right)\right)\right).
\tag{4}
\end{align*}

The first equality unravels the description of $\mathfrak{P}^{(2)}$ as a $(2,[1])$-identity second-order path; the second equality follows from Proposition~\ref{PDCHDUId}; the third equality applies the mapping $\mathrm{ip}^{([1],X)@}$ introduced in Definition~\ref{DPTQIp}; finally, the last equality applies the monomorphic Curry-Howard mapping at a path class.

By Proposition~\ref{PDEmb}, we have that
$$
\eta^{(2,1)\sharp}_{s}\left[
\eta^{(1,0)\sharp}_{s}\left[
\mathrm{T}_{\Sigma}(X)
\right]\right]
=
\eta^{(2,0)\sharp}\left[
\mathrm{T}_{\Sigma}(X)
\right].
$$
Moreover, taking into account that $\mathrm{CH}^{(2)}_{s}(\mathfrak{P}^{(2)})$ is a term in $\eta^{(2,0)\sharp}[\mathrm{T}_{\Sigma}(X)]$, we have that 
$
\mathrm{CH}^{(1)}_{s}(
\mathrm{ip}^{(1,X)@}_{s}(
P))
\in\eta^{(1,0)\sharp}[
\mathrm{T}_{\Sigma}(X)
]_{s}.
$

By Corollary~\ref{CCHId} $\mathrm{ip}^{(1,X)@}_{s}(P)$ is a $(1,0)$-identity path in $\mathrm{Pth}_{\boldsymbol{\mathcal{A}},s}$. Thus, there exists a term $Q\in\mathrm{T}_{\Sigma}(X)_{s}$ such that
$$
\mathrm{ip}^{(1,X)@}_{s}\left(
P
\right)
=
\mathrm{ip}^{(1,0)\sharp}_{s}\left(
Q
\right).
$$

The following chain of equalities holds 
\allowdisplaybreaks
\begin{align*}
\mathrm{CH}^{(1)}_{s}\left(
\mathrm{ip}^{(1,X)@}_{s}\left(
P
\right)\right)
&=
\mathrm{CH}^{(1)}_{s}\left(\mathrm{ip}^{(1,0)\sharp}_{s}\left(Q
\right)\right)
\tag{1}
\\&=
\eta^{(1,0)\sharp}_{s}(Q).
\tag{2}
\end{align*}

The first equation recovers the description of $\mathrm{ip}^{(1,X)@}_{s}(
P)$ as a $(1,0)$-identity path presented above; finally, the last equation follows from Proposition~\ref{PCHId}.

Thus, by Lemma~\ref{LWCong} we have that 
$$
\left(P,
\eta^{(1,0)\sharp}_{s}\left(
Q
\right)\right)
\in 
\Theta^{[1]}_{s}.
$$

Therefore, the following chain of equalities holds
\allowdisplaybreaks
\begin{align*}
\mathfrak{P}^{(2)}&=
\mathrm{ip}^{(2,[1])\sharp}_{s}\left(
\left[
\eta^{(1,0)\sharp}_{s}(Q)
\right]_{s}
\right)
\tag{1}
\\&=
\mathrm{ip}^{(2,[1])\sharp}_{s}\left(
\left[\mathrm{CH}^{(1)}_{s}\left(
\mathrm{ip}^{(1,0)\sharp}_{s}\left(
Q\right)\right)
\right]_{s}
\right)
\tag{2}
\\&=
\mathrm{ip}^{(2,[1])\sharp}_{s}\left(
\mathrm{CH}^{[1]}_{s}\left(
\left[\mathrm{ip}^{(1,0)\sharp}_{s}\left(
Q
\right)\right]_{s}
\right)\right)
\tag{3}
\\&=
\mathrm{ip}^{(2,[1])\sharp}_{s}\left(
\mathrm{CH}^{[1]}_{s}\left(
\mathrm{ip}^{([1],0)\sharp}_{s}\left(
Q
\right)\right)\right)
\tag{4}
\\&=
\mathrm{ip}^{(2,0)\sharp}_{s}(Q).
\tag{5}
\end{align*}

The first equality follows from Lemma~\ref{LWCong}; the second equality follows from Proposition~\ref{PCHId}; the third equality recovers the mapping $\mathrm{CH}^{[1]}$ introduced in Definition~\ref{DPTQCH}; the fourth equality recovers the mapping $\mathrm{ip}^{([1],0)\sharp}$ introduced in Definition~\ref{DCHUZ}; finally, the last equality recovers the mapping $\mathrm{ip}^{(2,0)\sharp}$ introduced in Definition~\ref{DDScTgZ}.

Therefore we can affirm that (i) $\mathfrak{P}^{(2)}$ is a $(2,0)$-identity second-order path.

Replacing the $(0,2)$-source by the $(0,2)$-target in the previous development follows easily since they coincide on $(2,0)$-identity second-order paths.

This finishes the proof.
\end{proof}

We will prove below that the non-$(2,[1])$-identity second-order paths containing   second-order echelons also have a nice characterization through the second-order Curry-Howard mapping. But, to do so, we first need to define some distinguished families of $S$-sorted subsets of terms in $\mathrm{T}_{\Sigma^{\boldsymbol{\mathcal{A}}^{(2)}}}(X)$ that are related to the binary operation of $1$-composition and to the $S$-sorted set of second-order rewrite rules.

\begin{definition}\label{DDA}
For the free $\Sigma^{\boldsymbol{\mathcal{A}}^{(2)}}$-algebra $\mathbf{T}_{\Sigma^{\boldsymbol{\mathcal{A}}^{(2)}}}(X)$
we define the following $S$-sorted subsets.

\textsf{(i)} The \emph{initial image of $\mathcal{A}^{(2)}$}, denoted by $\eta^{(2,\mathcal{A}^{(2)})}[\mathcal{A}^{(2)}]^{\mathrm{int}}$, is the $S$-sorted subset of $\mathrm{T}_{\Sigma^{\boldsymbol{\mathcal{A}}^{(2)}}}(X)$ defined, for every sort $s\in S$, as
$$
\eta^{(2,\mathcal{A}^{(2)})}[\mathcal{A}^{(2)}]^{\mathrm{int}}_{s}=
\left\lbrace
Q
\circ^{1\mathbf{T}_{\Sigma^{\boldsymbol{\mathcal{A}}^{(2)}}}(X)}_{s}
\mathfrak{p}^{(2)\mathbf{T}_{\Sigma^{\boldsymbol{\mathcal{A}}^{(2)}}}(X)}
\,\middle|\,
\begin{gathered}
Q\in\mathrm{T}_{\Sigma^{\boldsymbol{\mathcal{A}}^{(2)}}}(X)_{s},
\\
\mathfrak{p}^{(2)}\in\mathcal{A}^{(2)}_{s}.
\end{gathered}
\right\rbrace.
$$

\textsf{(ii)} The \emph{prescriptive image of $\mathcal{A}^{(2)}$}, denoted by $\eta^{(2,\mathcal{A}^{(2)})}[\mathcal{A}^{(2)}]^{\mathrm{pct}}$, is the $S$-sorted subset of $\mathrm{T}_{\Sigma^{\boldsymbol{\mathcal{A}}^{(2)}}}(X)$ defined, for every sort $s\in S$, as
\begin{multline*}
\eta^{(2,\mathcal{A}^{(2)})}[\mathcal{A}^{(2)}]^{\mathrm{pct}}_{s}\\=
\left\lbrace
\begin{gathered}
\mathfrak{p}^{(2)\mathbf{T}_{\Sigma^{\boldsymbol{\mathcal{A}}^{(2)}}}(X)},
\\
Q
\circ^{1\mathbf{T}_{\Sigma^{\boldsymbol{\mathcal{A}}^{(2)}}}(X)}_{s}
\mathfrak{p}^{(2)\mathbf{T}_{\Sigma^{\boldsymbol{\mathcal{A}}^{(2)}}}(X)},
\\
\left(R
\circ^{1\mathbf{T}_{\Sigma^{\boldsymbol{\mathcal{A}}^{(2)}}}(X)}_{s}
\mathfrak{p}^{(2)\mathbf{T}_{\Sigma^{\boldsymbol{\mathcal{A}}^{(2)}}}(X)}
\right)
\circ^{1\mathbf{T}_{\Sigma^{\boldsymbol{\mathcal{A}}^{(2)}}}(X)}_{s}
Q,
\\
\mathfrak{p}^{(2)\mathbf{T}_{\Sigma^{\boldsymbol{\mathcal{A}}^{(2)}}}(X)}
\circ^{1\mathbf{T}_{\Sigma^{\boldsymbol{\mathcal{A}}^{(2)}}}(X)}_{s}
Q
\end{gathered}
\,
\middle|
\,
\begin{gathered}
R,Q\in\mathrm{T}_{\Sigma^{\boldsymbol{\mathcal{A}}^{(2)}}}(X)_{s},
\\
\mathfrak{p}^{(2)}\in\mathcal{A}^{(2)}_{s}
\end{gathered}
\right\rbrace.
\end{multline*}

\textsf{(iii)} The \emph{non-initial image of $\mathcal{A}^{(2)}$}, denoted by $\eta^{(2,\mathcal{A}^{(2)})}[\mathcal{A}^{(2)}]^{\neg\mathrm{int}}$, is the $S$-sorted subset of $\mathrm{T}_{\Sigma^{\boldsymbol{\mathcal{A}}^{(2)}}}(X)$ defined, for every sort $s\in S$, as
\begin{multline*}
\eta^{(2,\mathcal{A}^{(2)})}[\mathcal{A}^{(2)}]^{\neg\mathrm{int}}_{s}\\=
\left\lbrace
\begin{gathered}
\left(R
\circ^{1\mathbf{T}_{\Sigma^{\boldsymbol{\mathcal{A}}^{(2)}}}(X)}_{s}
\mathfrak{p}^{(2)\mathbf{T}_{\Sigma^{\boldsymbol{\mathcal{A}}^{(2)}}}(X)}
\right)
\circ^{1\mathbf{T}_{\Sigma^{\boldsymbol{\mathcal{A}}^{(2)}}}(X)}_{s}
Q,
\\
\mathfrak{p}^{(2)\mathbf{T}_{\Sigma^{\boldsymbol{\mathcal{A}}^{(2)}}}(X)}
\circ^{1\mathbf{T}_{\Sigma^{\boldsymbol{\mathcal{A}}^{(2)}}}(X)}_{s}
Q
\end{gathered}
\,
\middle|
\,
\right.
\\
\left.
\begin{gathered}
R\in\mathrm{T}_{\Sigma^{\boldsymbol{\mathcal{A}}^{(2)}}}(X)_{s},
\\
Q\in \mathrm{T}_{\Sigma^{\boldsymbol{\mathcal{A}}^{(2)}}}(X)_{s}-\eta^{(2,\mathcal{A}^{(2)})}[\mathcal{A}^{(2)}]^{\mathrm{pct}}_{s},
\\ 
\mathfrak{p}^{(2)}\in\mathcal{A}^{(2)}_{s}
\end{gathered}
\right\rbrace.
\end{multline*}

\end{definition}

To prove the following lemmas one, basically, only needs to unravel the definition of the second-order Curry-Howard mapping at the suitable second-order path. Thus, we leave the proof to the reader.

\begin{lemma}\label{LDCHEch}
Let $s$ be a sort in $S$ and $\mathfrak{P}^{(2)}$ a second-order path in $\mathrm{Pth}_{\boldsymbol{\mathcal{A}}^{(2)},s}$. The following statements are equivalent.
\begin{enumerate}
\item[(i)] $\mathfrak{P}^{(2)}$ contains   second-order echelons.
\item[(ii)] $\mathrm{CH}^{(2)}_{s}(\mathfrak{P}^{(2)})\in \eta^{(2,\mathcal{A}^{(2)})}[\mathcal{A}^{(2)}]^{\mathrm{pct}}_{s}$.
\end{enumerate}
\end{lemma}

\begin{lemma}\label{LDCHDEch}
Let $s$ be a sort in $S$ and $\mathfrak{P}^{(2)}$ a second-order path in $\mathrm{Pth}_{\boldsymbol{\mathcal{A}}^{(2)},s}$. The following statements are equivalent.
\begin{enumerate}
\item[(i)] $\mathfrak{P}^{(2)}$ is a   second-order echelon.
\item[(ii)] $\mathrm{CH}^{(2)}_{s}(\mathfrak{P}^{(2)})\in \eta^{(2,\mathcal{A}^{(2)})}[\mathcal{A}^{(2)}]_{s}$.
\end{enumerate}
\end{lemma}

\begin{lemma}\label{LDCHEchInt}
Let $s$ be a sort in $S$ and $\mathfrak{P}^{(2)}$ a second-order path in $\mathrm{Pth}_{\boldsymbol{\mathcal{A}}^{(2)},s}$. The following statements are equivalent.
\begin{enumerate}
\item[(i)] $\bb{\mathfrak{P}^{(2)}}>1$ and $\mathfrak{P}^{(2)}$ has its first second-order echelon on its first step.
\item[(ii)] $\mathrm{CH}^{(2)}_{s}(\mathfrak{P}^{(2)})\in \eta^{(2,\mathcal{A}^{(2)})}[\mathcal{A}^{(2)}]^{\mathrm{int}}_{s}$.
\end{enumerate}
\end{lemma}

\begin{lemma}\label{LDCHEchNInt}
Let $s$ be a sort in $S$ and $\mathfrak{P}^{(2)}$ a second-order path in $\mathrm{Pth}_{\boldsymbol{\mathcal{A}}^{(2)},s}$. The following statements are equivalent.
\begin{enumerate}
\item[(i)] $\bb{\mathfrak{P}^{(2)}}>1$ and $\mathfrak{P}^{(2)}$ has its first second-order echelon on a step different from the initial one.
\item[(ii)] $\mathrm{CH}^{(2)}_{s}(\mathfrak{P}^{(2)})\in \eta^{(2,\mathcal{A}^{(2)})}[\mathcal{A}^{(2)}]^{\neg\mathrm{int}}_{s}$.
\end{enumerate}
\end{lemma}

We will prove below that the echelonless second-order paths also have a nice characterization through the second-order Curry-Howard mapping. But, to do so, we first need to define some distinguished families of $S$-sorted subsets of terms in $\mathrm{T}_{\Sigma^{\boldsymbol{\mathcal{A}}^{(2)}}}(X)$ that are related to the binary operation of $1$-composition and to the operations in the categorial signature determined by $\mathcal{A}$.

\begin{definition}\label{DDAII}
For the free $\Sigma^{\boldsymbol{\mathcal{A}}^{(2)}}$-algebra $\mathbf{T}_{\Sigma^{\boldsymbol{\mathcal{A}}^{(2)}}}(X)$
we define the following $S$-sorted subsets.

\textsf{(i)} Let $s$ be a sort in $S$, $\mathbf{s}$ a word in $S^{\star}-\{\lambda\}$, and $\tau\in\Sigma^{\boldsymbol{\mathcal{A}}^{(2)}}_{\mathbf{s}, s}$. Then we denote by $\mathcal{T}(\tau,\mathrm{T}_{\Sigma^{\boldsymbol{\mathcal{A}}^{(2)}}}(X)) = (\mathcal{T}(\tau,\mathrm{T}_{\Sigma^{\boldsymbol{\mathcal{A}}^{(2)}}}(X))_{n})_{n\in\mathbb{N}-\{0\}}$ the family of subsets of $\mathrm{T}_{\Sigma^{\boldsymbol{\mathcal{A}}^{(2)}}}(X)_{s}$ recursively defined as:
\[
\mathcal{T}\left(
\tau,\mathrm{T}_{\Sigma^{\boldsymbol{\mathcal{A}}^{(2)}}}(X)
\right)_{1}
=
\left\lbrace\tau^{\mathbf{T}_{\Sigma^{\boldsymbol{\mathcal{A}}^{(2)}}}(X)}
\left(\left(P_{j}
\right)_{j\in\bb{\mathbf{s}}}
\right)
\, 
\middle|
\,
\begin{gathered}
(P_{j})_{j\in\bb{\mathbf{s}}}
\in\mathrm{T}_{\Sigma^{\boldsymbol{\mathcal{A}}^{(2)}}}(X)_{\mathbf{s}},
\\
\exists j\in\bb{\mathbf{s}}\,
(P_{j}\not\in \eta^{(2,1)\sharp}
[\mathrm{PT}_{\boldsymbol{\mathcal{A}}}]_{s_{j}})
\end{gathered}
\right\rbrace,
\]
\begin{multline*}
\mathcal{T}\left(
\tau,\mathrm{T}_{\Sigma^{\boldsymbol{\mathcal{A}}^{(2)}}}(X)
\right)_{n+1}
\\=
\left\lbrace
Q\circ_{s}^{1\mathbf{T}_{\Sigma^{\boldsymbol{\mathcal{A}}^{(2)}}}(X)}
\tau^{\mathbf{T}_{\Sigma^{\boldsymbol{\mathcal{A}}^{(2)}}}(X)}
\left(\left(P_{j}
\right)_{j\in\bb{\mathbf{s}}}\right)
\, 
\middle|
\,
\begin{gathered}
Q\in\mathcal{T}(\tau,\mathrm{T}_{\Sigma^{\boldsymbol{\mathcal{A}}^{(2)}}}(X))_{n},\\
(P_{j})_{j\in\bb{\mathbf{s}}}
\in\mathrm{T}_{\Sigma^{\boldsymbol{\mathcal{A}}^{(2)}}}(X)_{\mathbf{s}},
\\
\exists j\in\bb{\mathbf{s}}\,
(P_{j}\not\in\eta^{(2,1)\sharp}
[\mathrm{PT}_{\boldsymbol{\mathcal{A}}}]_{s_{j}})
\end{gathered}
\right\rbrace.
\end{multline*}

\textsf{(ii)} Let $s$ be a sort in $S$, $\mathbf{s}$ a word in $S^{\star}-\{\lambda\}$, and $\tau\in\Sigma^{\boldsymbol{\mathcal{A}}}_{\mathbf{s}, s}$. Then we let $\mathcal{T}(\tau,\mathrm{T}_{\Sigma^{\boldsymbol{\mathcal{A}}^{(2)}}}(X))^{\star}$ and $\mathcal{T}(\tau,\mathrm{T}_{\Sigma^{\boldsymbol{\mathcal{A}}^{(2)}}}(X))^{+}$ stand for the subsets of $\mathrm{T}_{\Sigma^{\boldsymbol{\mathcal{A}}}}(X)_{s}$ given, respectively, by
\begin{align*}
\mathcal{T}\left(
\tau,\mathrm{T}_{\Sigma^{\boldsymbol{\mathcal{A}}^{(2)}}}(X)
\right)^{\star}
&=
\bigcup_{n\in\mathbb{N}-\{0\}}
\mathcal{T}\left(
\tau,\mathrm{T}_{\Sigma^{\boldsymbol{\mathcal{A}}^{(2)}}}(X)
\right)_{n},\,\,\text{and}
\\
\mathcal{T}\left(
\tau,\mathrm{T}_{\Sigma^{\boldsymbol{\mathcal{A}}^{(2)}}}(X)
\right)^{+}
&=
\bigcup_{n\in\mathbb{N}-\{0,1\}}
\mathcal{T}\left(
\tau,\mathrm{T}_{\Sigma^{\boldsymbol{\mathcal{A}}^{(2)}}}(X)
\right)_{n}.
\end{align*}

\textsf{(iii)} The \emph{head-constant and coherent $S$-sorted subset of $\mathrm{T}_{\Sigma^{\boldsymbol{\mathcal{A}}^{(2)}}}(X)$}, written  $\mathrm{T}_{\Sigma^{\boldsymbol{\mathcal{A}}^{(2)}}}(X)^{
\mathsf{HdC}\And\mathsf{C}
}$, is  defined, for every sort $s\in S$, as
$$
\mathrm{T}_{\Sigma^{\boldsymbol{\mathcal{A}}^{(2)}}}(X)
^{\mathsf{HdC}\And\mathsf{C}}_{s}=
\bigcup_{\tau\in\Sigma^{\boldsymbol{\mathcal{A}}}_{\neq\lambda, s}}
\mathcal{T}
\left(
\tau,\mathrm{T}_{\Sigma^{\boldsymbol{\mathcal{A}}^{(2)}}}(X)
\right)_{1}.
$$

\textsf{(iv)} The \emph{head-constant and non-coherent $S$-sorted subset of $\mathrm{T}_{\Sigma^{\boldsymbol{\mathcal{A}}^{(2)}}}(X)$}, written $\mathrm{T}_{\Sigma^{\boldsymbol{\mathcal{A}}^{(2)}}}(X)^{\mathsf{HdC}\And\neg\mathsf{C}}$, is  defined, for every sort $s\in S$, as
$$
\mathrm{T}_{\Sigma^{\boldsymbol{\mathcal{A}}^{(2)}}}(X)^{\mathsf{HdC}\And\neg\mathsf{C}}_{s}=
\bigcup_{\tau\in\Sigma^{\boldsymbol{\mathcal{A}}}_{\neq\lambda, s}}
\mathcal{T}
\left(
\tau,\mathrm{T}_{\Sigma^{\boldsymbol{\mathcal{A}}^{(2)}}}(X)
\right)^{+}.
$$

\textsf{(v)} The \emph{head-constant $S$-sorted subset of $\mathrm{T}_{\Sigma^{\boldsymbol{\mathcal{A}}^{(2)}}}(X)$}, written  $\mathrm{T}_{\Sigma^{\boldsymbol{\mathcal{A}}^{(2)}}}(X)^{\mathsf{HdC}}$, is  defined, for every sort $s\in S$, as
$$
\mathrm{T}_{\Sigma^{\boldsymbol{\mathcal{A}}^{(2)}}}(X)^{\mathsf{HdC}}_{s}=
\bigcup_{\tau\in\Sigma^{\boldsymbol{\mathcal{A}}}_{\neq\lambda, s}}
\mathcal{T}
\left(
\tau,\mathrm{T}_{\Sigma^{\boldsymbol{\mathcal{A}}^{(2)}}}(X)
\right)^{\star}.
$$

\textsf{(iv)} The \emph{non head-constant $S$-sorted subset of $\mathrm{T}_{\Sigma^{\boldsymbol{\mathcal{A}}^{(2)}}}(X)$}, written  $\mathrm{T}_{\Sigma^{\boldsymbol{\mathcal{A}}^{(2)}}}(X)^{\neg\mathsf{HdC}}$, is  defined, for every sort $s\in S$, as
\begin{multline*}
\mathrm{T}_{\Sigma^{\boldsymbol{\mathcal{A}}^{(2)}}}(X)^{\neg\mathsf{HdC}}_{s}
\\=
\left\lbrace
Q\circ^{1\mathbf{T}_{\Sigma^{\boldsymbol{\mathcal{A}}^{(2)}}}(X)}
\tau^{\mathbf{T}_{\Sigma^{\boldsymbol{\mathcal{A}}^{(2)}}}(X)}
\left(\left(
P_{j}
\right)_{j\in\bb{\mathbf{s}}}
\right)
\, 
\middle|
\,
\begin{gathered}
\mathbf{s}\in S^{\star}-\{\lambda\}, \tau\in\Sigma^{\boldsymbol{\mathcal{A}}}_{\mathbf{s},s},\\
(P_{j})_{j\in\bb{\mathbf{s}}}\in\mathrm{T}_{\Sigma^{\boldsymbol{\mathcal{A}}^{(2)}}}(X)_{\mathbf{s}},\\
\exists j\in\bb{\mathbf{s}}\,(P_{j}\not \in\eta^{(2,1)\sharp}[\mathrm{PT}_{\boldsymbol{\mathcal{A}}}]_{s_{j}}),\\
Q\in\mathrm{T}_{\Sigma^{\boldsymbol{\mathcal{A}}^{(2)}}}(X)_{s},\\
Q\not\in\eta^{(2,1)\sharp}[\mathrm{PT}_{\boldsymbol{\mathcal{A}}}]_{s},\\
Q\not\in\eta^{(2,\mathcal{A}^{(2)})}[\mathcal{A}^{(2)}]^{\mathrm{pct}}_{s},\\
Q\not\in\mathcal{T}(\tau,\mathrm{T}_{\Sigma^{\boldsymbol{\mathcal{A}}^{(2)}}}(X))^{\star}.
\end{gathered}
\right\rbrace.
\end{multline*}
\end{definition}

To prove the following lemmas one, basically, only needs to unravel the definition of the second-order Curry-Howard mapping at the suitable second-order path. Thus, we leave the proof to the reader.

\begin{lemma}\label{LDCHNEch}
Let $s$ be a sort in $S$ and $\mathfrak{P}^{(2)}$ a second-order path in $\mathrm{Pth}_{\boldsymbol{\mathcal{A}}^{(2)},s}$. The following statements are equivalent.
\begin{enumerate}
\item[(i)] $\mathfrak{P}^{(2)}$ is an  echelonless second-order path.
\item[(ii)] $\mathrm{CH}^{(2)}_{s}(\mathfrak{P}^{(2)})\in \mathrm{T}_{\Sigma^{\boldsymbol{\mathcal{A}}^{(2)}}}(X)_{s}-(\eta^{(2,1)\sharp}[\mathrm{PT}_{\boldsymbol{\mathcal{A}}}]_{s}\cup
\eta^{(2,\mathcal{A}^{(2)})}[\mathcal{A}^{(2)}]^{\mathrm{pct}}_{s}
)
$.
\end{enumerate}
\end{lemma}

\begin{lemma}\label{LDCHNEchNHd}
Let $s$ be a sort in $S$ and $\mathfrak{P}^{(2)}$ a second-order path in $\mathrm{Pth}_{\boldsymbol{\mathcal{A}}^{(2)},s}$.The following statements are equivalent.
\begin{enumerate}
\item[(i)] $\mathfrak{P}^{(2)}$ is an echelonless second-order path that is not head-constant.
\item[(ii)] $\mathrm{CH}^{(2)}_{s}(\mathfrak{P}^{(2)})\in \mathrm{T}_{\Sigma^{\boldsymbol{\mathcal{A}}^{(2)}}}(X)^{\neg\mathsf{HdC}}_{s}
$.
\end{enumerate}
\end{lemma}

\begin{lemma}\label{LDCHNEchHd}
Let $s$ be a sort in $S$ and $\mathfrak{P}^{(2)}$ a second-order path in $\mathrm{Pth}_{\boldsymbol{\mathcal{A}}^{(2)},s}$. The following statements are equivalent.
\begin{enumerate}
\item[(i)] $\mathfrak{P}^{(2)}$  is a head-constant echelonless second-order path.
\item[(ii)] $\mathrm{CH}^{(2)}_{s}(\mathfrak{P}^{(2)})\in \mathrm{T}_{\Sigma^{\boldsymbol{\mathcal{A}}^{(2)}}}(X)^{\mathsf{HdC}}_{s}
$.
\end{enumerate}
\end{lemma}

\begin{lemma}\label{LDCHNEchHdNC}
Let $s$ be a sort in $S$ and $\mathfrak{P}^{(2)}$ a second-order path in $\mathrm{Pth}_{\boldsymbol{\mathcal{A}}^{(2)},s}$. The following statements are equivalent.
\begin{enumerate}
\item[(i)] $\mathfrak{P}^{(2)}$ is a head-constant echelonless second-order path that is not coherent.
\item[(ii)] $\mathrm{CH}^{(2)}_{s}(\mathfrak{P}^{(2)})\in \mathrm{T}_{\Sigma^{\boldsymbol{\mathcal{A}}^{(2)}}}(X)^{\mathsf{HdC}\And\neg\mathsf{C}}_{s}
$.
\end{enumerate}
\end{lemma}

\begin{lemma}\label{LDCHNEchHdC}
Let $s$ be a sort in $S$ and $\mathfrak{P}^{(2)}$ a second-order path in $\mathrm{Pth}_{\boldsymbol{\mathcal{A}}^{(2)},s}$. The following statements are equivalent.
\begin{enumerate}
\item[(i)] $\mathfrak{P}^{(2)}$ is a  coherent head-constant echelonless second-order path.
\item[(ii)] $\mathrm{CH}^{(2)}_{s}(\mathfrak{P}^{(2)})\in \mathrm{T}_{\Sigma^{\boldsymbol{\mathcal{A}}^{(2)}}}(X)^{\mathsf{HdC}\And\mathsf{C}}_{s}
$.
\end{enumerate}
\end{lemma}

Another interesting result implied by the just stated lemmas is that the second-order Curry-Howard mapping restricted to the subset of $\mathrm{Pth}_{\boldsymbol{\mathcal{A}}^{(2)}}$ formed by the second-order echelons is an injective mapping.

\begin{restatable}{proposition}{PDCHEch}
\label{PDCHEch} Let $s$ be a sort in $S$and $(\mathfrak{P}^{(2)},\mathfrak{Q}^{(2)})\in\mathrm{Ker}(\mathrm{CH}^{(2)})_{s}$. If $\mathfrak{P}^{(2)}$ or $\mathfrak{Q}^{(2)}$ is a second-order echelon, then $\mathfrak{Q}^{(2)}=\mathfrak{P}^{(2)}$.
\end{restatable}
\begin{proof}
Assume that $\mathfrak{P}^{(2)}$ is a second-order echelon. That is $\mathfrak{P}^{(2)}=\mathrm{ech}^{(2,\mathcal{A}^{(2)})}_{s}(\mathfrak{p}^{(2)})$, for some second-order rewrite rule $\mathfrak{p}^{(2)}\in\mathcal{A}^{(2)}_{s}$. Hence, following Definition~\ref{DDCH}, 
\[
\mathrm{CH}^{(2)}_{s}\left(
\mathfrak{P}^{(2)}
\right)
=
\mathfrak{p}^{(2)\mathbf{T}_{\Sigma^{\boldsymbol{\mathcal{A}}^{(2)}}}(X)}.
\]

Since $(\mathfrak{P}^{(2)},\mathfrak{Q}^{(2)})\in\mathrm{Ker}(\mathrm{CH}^{(2)})_{s}$, we have that $\mathrm{CH}^{(2)}_{s}(\mathfrak{Q}^{(2)})$ is a term in $\eta^{(2,\mathcal{A}^{(2)})}[\mathcal{A}^{(2)}]_{s}$. Hence, by Lemma~\ref{LDCHDEch}, $\mathfrak{Q}^{(2)}$ is a second-order echelon. Since $\mathrm{CH}^{(2)}_{s}(\mathfrak{Q}^{(2)})=\mathrm{CH}^{(2)}_{s}(\mathfrak{P}^{(2)})$, we can affirm that the unique second-order rewrite rule appearing in both $\mathfrak{P}^{(2)}$ and $\mathfrak{Q}^{(2)}$ is the same. Consequently $\mathfrak{P}^{(2)}=\mathfrak{Q}^{(2)}$.
\end{proof}

Before proving that $\mathrm{Ker}(\mathrm{CH}^{(2)})$ is a $\Sigma^{\boldsymbol{\mathcal{A}}^{(2)}}$-congruence on $\mathbf{Pth}_{\boldsymbol{\mathcal{A}}^{(2)}}$, we need to prove, in the following lemma, that  whenever two second-order paths of the same sort are in the kernel of the second-order Curry-Howard mapping, then they have the same length, the same $([1],2)$-source and the same $([1],2)$-target.

\begin{restatable}{lemma}{LDCH}
\label{LDCH}
Let $s$ be a sort in $S$ and $(\mathfrak{P}^{(2)},\mathfrak{Q}^{(2)})\in\mathrm{Ker}(\mathrm{CH}^{(2)})_{s}$. The following statements hold.
\begin{enumerate}
\item[(i)] $\bb{\mathfrak{P}^{(2)}}=\bb{\mathfrak{Q}^{(2)}}$;
\item[(ii)] $\mathrm{sc}^{([1],2)}_{s}(\mathfrak{P}^{(2)})=\mathrm{sc}^{([1],2)}_{s}(\mathfrak{Q}^{(2)})$;
\item[(iii)] $\mathrm{tg}^{([1],2)}_{s}(\mathfrak{P}^{(2)})=\mathrm{tg}^{([1],2)}_{s}(\mathfrak{Q}^{(2)})$.
\end{enumerate}
\end{restatable}
\begin{proof}
If either $\mathfrak{P}^{(2)}$ or $\mathfrak{Q}^{(2)}$ is a $(2,[1])$-identity second-order path, then the statement follows according to Corollary~\ref{CDCHUId}. Therefore, we can assume that neither $\mathfrak{P}^{(2)}$ nor $\mathfrak{Q}^{(2)}$ are $(2,[1])$-identity second-order paths.

We prove the statements by Artinian induction on $(\coprod\mathrm{Pth}_{\boldsymbol{\mathcal{A}}^{(2)}}, \leq_{\mathbf{Pth}_{\boldsymbol{\mathcal{A}}^{(2)}}})$.

\textsf{Base step of the Artinian induction.}

Let $(\mathfrak{P}^{(2)}, s)$ be a minimal element of $(\coprod\mathrm{Pth}_{\boldsymbol{\mathcal{A}}^{(2)}}, \leq_{\mathbf{Pth}_{\boldsymbol{\mathcal{A}}^{(2)}}})$. Then, by Proposition~\ref{PDMinimal} and taking into account that we are assuming that $\mathfrak{P}^{(2)}$ is a non-$(2,[1])$-identity second-order path, we have that $\mathfrak{P}^{(2)}$ is a   second-order echelon. The statement follows in virtue of Proposition~\ref{PDCHEch}.

This completes the base step.

\textsf{Inductive step of the Artinian induction.}

Let $(\mathfrak{P}^{(2)},s)$ be a non-minimal element of $(\coprod
\mathrm{Pth}_{\boldsymbol{\mathcal{A}}^{(2)}},
\leq_{\mathbf{Pth}_{\boldsymbol{\mathcal{A}}^{(2)}}})
$. Let us suppose that, for every sort $t\in S$ and every second-order path $\mathfrak{P}'^{(2)}$ in $\mathrm{Pth}_{\boldsymbol{\mathcal{A}}^{(2)},t}$, if $(\mathfrak{P}'^{(2)},t)<_{\mathbf{Pth}_{\boldsymbol{\mathcal{A}}^{(2)}}}
(\mathfrak{P}^{(2)},s)
$, then the statement holds for $\mathfrak{P}'^{(2)}$, i.e., for every second-order path $\mathfrak{Q}'^{(2)}$ in $\mathrm{Pth}_{\boldsymbol{\mathcal{A}}^{(2)},t}$, if $(\mathfrak{P}'^{(2)}, \mathfrak{Q}'^{(2)})\in\mathrm{Ker}(\mathrm{CH}^{(2)})_{t}$, then $\mathfrak{P}'^{(2)}$ and $\mathfrak{Q}'^{(2)}$ have the same length, the same $([1],2)$-source and the same $([1],2)$-target.

Since $(\mathfrak{P}^{(2)},s)$ is a non-minimal element of $(\coprod\mathrm{Pth}_{\boldsymbol{\mathcal{A}}^{(2)}}, \leq_{\mathbf{Pth}_{\boldsymbol{\mathcal{A}}^{(2)}}})$ and taking into account that $\mathfrak{P}^{(2)}$ is a non-$(2,[1])$-identity second-order path, we have, by Lemma~\ref{LDOrdI}, that $\mathfrak{P}^{(2)}$ is either (1)~a second-order path of length strictly greater than one containing at least one   second-order echelon or (2)~an echelonless second-order path.

If~(1), then let $i\in\bb{\mathfrak{P}^{(2)}}$ be the first index for which the one-step subpath $\mathfrak{P}^{(2),i,i}$ is a   second-order echelon. We distinguish two cases accordingly.

If $i=0$, i.e., if $\mathfrak{P}^{(2)}$ has its first   second-order echelon on its first step, then, according to Definition~\ref{DDCH}, we have that
$$
\mathrm{CH}^{(2)}_{s}\left(
\mathfrak{P}^{(2)}
\right)=
\mathrm{CH}^{(2)}_{s}\left(
\mathfrak{P}^{(2),1,\bb{\mathfrak{P}^{(2)}}-1}
\right)
\circ_{s}^{1\mathbf{T}_{\Sigma^{\boldsymbol{\mathcal{A}}^{(2)}}}(X)}
\mathrm{CH}^{(2)}_{s}\left(
\mathfrak{P}^{(2),0,0}
\right).
$$

Since $\mathrm{CH}^{(2)}_{s}(\mathfrak{P}^{(2)})\in\eta^{(2,\mathcal{A}^{(2)})}[\mathcal{A}^{(2)}]^{\mathrm{int}}_{s}$ and $(\mathfrak{P}^{(2)},\mathfrak{Q}^{(2)})\in\mathrm{Ker}(\mathrm{CH}^{(2)})_{s}$, we have, by Lemma~\ref{LDCHEchInt}, that $\mathfrak{Q}^{(2)}$ is a path of length strictly greater than one containing its first   second-order echelon on its first step. 

Thus, according to Definition~\ref{DDCH}, we have that
$$
\mathrm{CH}^{(2)}_{s}\left(
\mathfrak{Q}^{(2)}
\right)
=
\mathrm{CH}^{(2)}_{s}\left(
\mathfrak{Q}^{(2),1,\bb{\mathfrak{Q}^{(2)}}-1}
\right)
\circ_{s}^{1\mathbf{T}_{\Sigma^{\boldsymbol{\mathcal{A}}^{(2)}}}(X)}
\mathrm{CH}^{(2)}_{s}\left(
\mathfrak{Q}^{(2),0,0}
\right).
$$

Since $(\mathfrak{P}^{(2)},\mathfrak{Q}^{(2)})\in\mathrm{Ker}(\mathrm{CH}^{(2)})_{s}$, we have that $(\mathfrak{P}^{(2),1,\bb{\mathfrak{P}^{(2)}}-1},\mathfrak{Q}^{(2),1,\bb{\mathfrak{Q}^{(2)}}-1})$ and $(\mathfrak{P}^{(2),0,0}, \mathfrak{Q}^{(2),0,0})$ are in $\mathrm{Ker}(\mathrm{CH}^{(2)})_{s}$. Note that, according to Definition~\ref{DDOrd}, we have that $(\mathfrak{P}^{(2),0,0},s)$ and $(\mathfrak{P}^{(2),1,\bb{\mathfrak{P}^{(2)}}-1},s)$ $\prec_{\mathbf{Pth}_{\boldsymbol{\mathcal{A}}^{(2)}}}$-precede $(\mathfrak{P}^{(2)},s)$. 

Therefore, by the inductive hypothesis, the second-order paths $\mathfrak{P}^{(2),0,0}$ and $\mathfrak{Q}^{(2),0,0}$, and the second-order paths $\mathfrak{P}^{(2),1,\bb{\mathfrak{P}^{(2)}}-1}$ and $\mathfrak{Q}^{(2),1,\bb{\mathfrak{Q}^{(2)}}-1}$ have, respectively, the same length, the same $([1],2)$-source and the same $([1],2)$-target.

Therefore, we have that
\leqnomode
\allowdisplaybreaks
\begin{align*}
\bb{\mathfrak{P}^{(2)}}&=
\bb{\mathfrak{P}^{(2),0,0}}+
\bb{\mathfrak{P}^{(2),1,\bb{\mathfrak{P}^{(2)}}-1}}
\tag{i}
\\&=
\bb{\mathfrak{Q}^{(2),0,0}}+
\bb{\mathfrak{Q}^{(2),1,\bb{\mathfrak{Q}^{(2)}}-1}}
\\&=\bb{\mathfrak{Q}^{(2)}};
\\
\mathrm{sc}^{([1],2)}_{s}\left(\mathfrak{P}^{(2)}
\right)
&=
\mathrm{sc}^{([1],2)}_{s}\left(
\mathfrak{P}^{(2),0,0}
\right)
\tag{ii}
\\&=
\mathrm{sc}^{([1],2)}_{s}\left(
\mathfrak{Q}^{(2),0,0}
\right)
\\&=
\mathrm{sc}^{([1],2)}_{s}\left(
\mathfrak{Q}^{(2)}
\right);
\\
\mathrm{tg}^{([1],2)}_{s}\left(
\mathfrak{P}^{(2)}\right)
&=
\mathrm{tg}^{([1],2)}_{s}\left(
\mathfrak{P}^{(2),1,\bb{\mathfrak{P}^{(2)}}-1}
\right)
\tag{iii}
\\&=
\mathrm{tg}^{([1],2)}_{s}\left(
\mathfrak{Q}^{(2),1,\bb{\mathfrak{Q}^{(2)}}-1}
\right)
\\&=
\mathrm{tg}^{([1],2)}_{s}\left(
\mathfrak{Q}^{(2)}\right).
\end{align*}

The case of $\mathfrak{P}^{(2)}$ being a second-order path of length strictly greater than one containing its first   second-order echelon on its first step follows.

If $i>0$, that is, if $\mathfrak{P}^{(2)}$ is a second-order path of length strictly greater than one containing its first   second-order echelon on a step different from the initial one, then, according to Definition~\ref{DDCH}, we have that
$$
\mathrm{CH}^{(2)}_{s}\left(
\mathfrak{P}^{(2)}
\right)=
\mathrm{CH}^{(2)}_{s}\left(
\mathfrak{P}^{(2),i,\bb{\mathfrak{P}^{(2)}}-1}
\right)
\circ_{s}^{1\mathbf{T}_{\Sigma^{\boldsymbol{\mathcal{A}}^{(2)}}}(X)}
\mathrm{CH}^{(2)}_{s}\left(
\mathfrak{P}^{(2),0,i-1}
\right).
$$

Since $\mathrm{CH}^{(2)}_{s}(\mathfrak{P}^{(2)})\in\eta^{(2,\mathcal{A}^{(2)})}[\mathcal{A}^{(2)}]^{\neg\mathrm{int}}_{s}$ and $(\mathfrak{P}^{(2)},\mathfrak{Q}^{(2)})\in\mathrm{Ker}(\mathrm{CH}^{(2)})_{s}$, we have, by Lemma~\ref{LDCHEchNInt}, that $\mathfrak{Q}^{(2)}$ is a second-order path of length strictly greater than one containing its first   second-order echelon on a step different from the initial one.

Thus, if $j\in\bb{\mathfrak{Q}^{(2)}}$ is the first index for which $\mathfrak{Q}^{(2),j,j}$ is a   second-order echelon, according to Definition~\ref{DDCH}, we have that
$$
\mathrm{CH}^{(2)}_{s}\left(
\mathfrak{Q}^{(2)}
\right)
=
\mathrm{CH}^{(2)}_{s}\left(
\mathfrak{Q}^{(2),j,\bb{\mathfrak{Q}^{(2)}}-1}
\right)
\circ_{s}^{1\mathbf{T}_{\Sigma^{\boldsymbol{\mathcal{A}}^{(2)}}}(X)}
\mathrm{CH}^{(2)}_{s}\left(
\mathfrak{Q}^{(2),0,j-1}
\right).
$$

Since $(\mathfrak{P}^{(2)},\mathfrak{Q}^{(2)})\in\mathrm{Ker}(\mathrm{CH}^{(2)})_{s}$, we have that $(\mathfrak{P}^{(2),i,\bb{\mathfrak{P}^{(2)}}-1},\mathfrak{Q}^{(2),j,\bb{\mathfrak{Q}^{(2)}}-1})$ and $(\mathfrak{P}^{(2),0,i-1}, \mathfrak{Q}^{(2),0,j-1})$ are in $\mathrm{Ker}(\mathrm{CH}^{(2)})_{s}$. Note that, according to Definition~\ref{DDOrd}, we have that $(\mathfrak{P}^{(2),0,i-1},s)$ and $(\mathfrak{P}^{(2),i,\bb{\mathfrak{P}^{(2)}}-1},s)$ $\prec_{\mathbf{Pth}_{\boldsymbol{\mathcal{A}}^{(2)}}}$-precede $(\mathfrak{P}^{(2)},s)$. 

Therefore, by the inductive hypothesis, the second-order paths $\mathfrak{P}^{(2),0,i-1}$ and $\mathfrak{Q}^{(2),0,j-1}$, and the second-order paths $\mathfrak{P}^{(2),i,\bb{\mathfrak{P}^{(2)}}-1}$ and $\mathfrak{Q}^{(2),j,\bb{\mathfrak{Q}^{(2)}}-1}$ have, respectively, the same length, the same $([1],2)$-source and the same $([1],2)$-target. In particular $i=j$.

Therefore, we have that
\leqnomode
\allowdisplaybreaks
\begin{align*}
\bb{\mathfrak{P}^{(2)}}&=
\bb{\mathfrak{P}^{(2),0,i-1}}+
\bb{\mathfrak{P}^{(2),i,\bb{\mathfrak{P}^{(2)}}-1}}
\tag{i}
\\&=
\bb{\mathfrak{Q}^{(2),0,j-1}}+
\bb{\mathfrak{Q}^{(2),j,\bb{\mathfrak{Q}^{(2)}}-1}}
\\&=\bb{\mathfrak{Q}^{(2)}};
\\
\mathrm{sc}^{([1],2)}_{s}\left(
\mathfrak{P}^{(2)}
\right)
&=
\mathrm{sc}^{([1],2)}_{s}\left(
\mathfrak{P}^{(2),0,i-1}
\right)
\tag{ii}
\\&=
\mathrm{sc}^{([1],2)}_{s}\left(
\mathfrak{Q}^{(2),0,j-1}
\right)
\\&=
\mathrm{sc}^{([1],2)}_{s}\left(\mathfrak{Q}^{(2)}
\right);
\\
\mathrm{tg}^{([1],2)}_{s}\left(
\mathfrak{P}^{(2)}
\right)
&=
\mathrm{tg}^{([1],2)}_{s}\left(
\mathfrak{P}^{(2),i,\bb{\mathfrak{P}^{(2)}}-1}
\right)
\tag{iii}
\\&=
\mathrm{tg}^{([1],2)}_{s}\left(
\mathfrak{Q}^{(2),j,\bb{\mathfrak{Q}^{(2)}}-1}
\right)
\\&=
\mathrm{tg}^{([1],2)}_{s}\left(
\mathfrak{Q}^{(2)}
\right)
.
\end{align*}

The case of $\mathfrak{P}^{(2)}$ being a second-order path of length strictly greater than one containing its first   second-order echelon on a step different from the initial one follows.

Case~(1) follows.

If~(2), i.e., if $\mathfrak{P}^{(2)}$ is an echelonless second-order path, it could be the case that (2.1) $\mathfrak{P}^{(2)}$ is an echelonless second-order path that is not head-constant, or (2.2) $\mathfrak{P}^{(2)}$ is a head-constant non-coherent echelonless second-order path or (2.3) $\mathfrak{P}^{(2)}$ is a head-constant coherent echelonless second-order path.

If~(2.1), let $i\in\bb{\mathfrak{P}^{(2)}}$ be the greatest index for which $\mathfrak{P}^{(2),0,i}$ is a head-constant second-order path. Then, according to Definition~\ref{DDCH}, we have that
$$
\mathrm{CH}^{(2)}_{s}\left(
\mathfrak{P}^{(2)}
\right)=
\mathrm{CH}^{(2)}_{s}\left(
\mathfrak{P}^{(2),i+1,\bb{\mathfrak{P}^{(2)}}-1}
\right)
\circ_{s}^{1\mathbf{T}_{\Sigma^{\boldsymbol{\mathcal{A}}^{(2)}}}(X)}
\mathrm{CH}^{(2)}_{s}\left(
\mathfrak{P}^{(2),0,i}
\right).
$$

Since $\mathrm{CH}^{(2)}_{s}(\mathfrak{P}^{(2)})\in\mathrm{T}_{\Sigma^{\boldsymbol{\mathcal{A}}^{(2)}}}(X)^{\neg\mathsf{HdC}}$ and $(\mathfrak{P}^{(2)},\mathfrak{Q}^{(2)})\in\mathrm{Ker}(\mathrm{CH}^{(2)})_{s}$, we have, by Lemma~\ref{LDCHNEchNHd}, that $\mathfrak{Q}^{(2)}$ is an echelonless second-order path  that is not head-constsnt.

Thus, if $j\in\bb{\mathfrak{Q}^{(2)}}$ is the greatest index for which $\mathfrak{Q}^{(2),0,j}$ is a head-constant second-order path, according to Definition~\ref{DDCH}, we have that
$$
\mathrm{CH}^{(2)}_{s}\left(
\mathfrak{Q}^{(2)}
\right)
=
\mathrm{CH}^{(2)}_{s}\left(
\mathfrak{Q}^{(2),j+1,\bb{\mathfrak{Q}^{(2)}}-1}
\right)
\circ_{s}^{1\mathbf{T}_{\Sigma^{\boldsymbol{\mathcal{A}}^{(2)}}}(X)}
\mathrm{CH}^{(2)}_{s}\left(
\mathfrak{Q}^{(2),0,j}
\right).
$$

Since $(\mathfrak{P}^{(2)},\mathfrak{Q}^{(2)})\in\mathrm{Ker}(\mathrm{CH}^{(2)})_{s}$, we have that $(\mathfrak{P}^{(2),i+1,\bb{\mathfrak{P}^{(2)}}-1},\mathfrak{Q}^{(2),j+1,\bb{\mathfrak{Q}^{(2)}}-1})$ and $(\mathfrak{P}^{(2),0,i}, \mathfrak{Q}^{(2),0,j})$ are in $\mathrm{Ker}(\mathrm{CH}^{(2)})_{s}$. Note that, according to Definition~\ref{DDOrd}, we have that $(\mathfrak{P}^{(2),0,i},s)$ and $(\mathfrak{P}^{(2),i+1,\bb{\mathfrak{P}^{(2)}}-1},s)$ $\prec_{\mathbf{Pth}_{\boldsymbol{\mathcal{A}}^{(2)}}}$-precede $(\mathfrak{P}^{(2)},s)$. 

Therefore, by the inductive hypothesis, the second-order paths $\mathfrak{P}^{(2),0,i}$ and $\mathfrak{Q}^{(2),0,j}$, and the second-order paths $\mathfrak{P}^{(2),i+1,\bb{\mathfrak{P}^{(2)}}-1}$ and $\mathfrak{Q}^{(2),j+1,\bb{\mathfrak{Q}^{(2)}}-1}$ have, respectively, the same length, the same $([1],2)$-source and the same $([1],2)$-target. In particular $i=j$.

Therefore, we have that
\leqnomode
\allowdisplaybreaks
\begin{align*}
\bb{\mathfrak{P}^{(2)}}&=
\bb{\mathfrak{P}^{(2),0,i}}+
\bb{\mathfrak{P}^{(2),i+1,\bb{\mathfrak{P}^{(2)}}-1}}
\tag{i}
\\&=
\bb{\mathfrak{Q}^{(2),0,j}}+
\bb{\mathfrak{Q}^{(2),j+1,\bb{\mathfrak{Q}^{(2)}}-1}}
\\&=\bb{\mathfrak{Q}^{(2)}};
\\
\mathrm{sc}^{([1],2)}_{s}\left(
\mathfrak{P}^{(2)}
\right)
&=
\mathrm{sc}^{([1],2)}_{s}\left(
\mathfrak{P}^{(2),0,i}
\right)
\tag{ii}
\\&=
\mathrm{sc}^{([1],2)}_{s}\left(
\mathfrak{Q}^{(2),0,j}
\right)
\\&=
\mathrm{sc}^{([1],2)}_{s}\left(
\mathfrak{Q}^{(2)}
\right)
;
\\
\mathrm{tg}^{([1],2)}_{s}\left(
\mathfrak{P}^{(2)}
\right)
&=
\mathrm{tg}^{([1],2)}_{s}\left(
\mathfrak{P}^{(2),i+1,\bb{\mathfrak{P}^{(2)}}-1}
\right)
\tag{iii}
\\&=
\mathrm{tg}^{([1],2)}_{s}\left(
\mathfrak{Q}^{(2),j+1,\bb{\mathfrak{Q}^{(2)}}-1}
\right)
\\&=
\mathrm{tg}^{([1],2)}_{s}\left(
\mathfrak{Q}^{(2)}
\right)
.
\end{align*}

The case of $\mathfrak{P}^{(2)}$ being an echelonless second-order path  that is not head-constant follows.

If~(2.2), let $i\in\bb{\mathfrak{P}^{(2)}}$ be the greatest index for which $\mathfrak{P}^{(2),0,i}$ is a coherent second-order path. Then, according to Definition~\ref{DDCH}, we have that
$$
\mathrm{CH}^{(2)}_{s}\left(
\mathfrak{P}^{(2)}
\right)=
\mathrm{CH}^{(2)}_{s}\left(
\mathfrak{P}^{(2),i+1,\bb{\mathfrak{P}^{(2)}}-1}
\right)
\circ_{s}^{1\mathbf{T}_{\Sigma^{\boldsymbol{\mathcal{A}}^{(2)}}}(X)}
\mathrm{CH}^{(2)}_{s}\left(
\mathfrak{P}^{(2),0,i}
\right).
$$

Since $\mathrm{CH}^{(2)}_{s}(\mathfrak{P}^{(2)})\in\mathrm{T}_{\Sigma^{\boldsymbol{\mathcal{A}}^{(2)}}}(X)^{\mathsf{HdC}\And\neg\mathsf{C}}$ and $(\mathfrak{P}^{(2)},\mathfrak{Q}^{(2)})\in\mathrm{Ker}(\mathrm{CH}^{(2)})_{s}$ we have, by Lemma~\ref{LDCHNEchHdNC} that $\mathfrak{Q}^{(2)}$ is an echelonless second-order path that is not head-constant.

Thus, if $j\in\bb{\mathfrak{Q}^{(2)}}$ is the greatest index for which $\mathfrak{Q}^{(2),0,j}$ is a coherent second-order path, then, according to Definition~\ref{DDCH}, we have that
$$
\mathrm{CH}^{(2)}_{s}\left(
\mathfrak{Q}^{(2)}
\right)
=
\mathrm{CH}^{(2)}_{s}\left(
\mathfrak{Q}^{(2),j+1,\bb{\mathfrak{Q}^{(2)}}-1}
\right)
\circ_{s}^{1\mathbf{T}_{\Sigma^{\boldsymbol{\mathcal{A}}^{(2)}}}(X)}
\mathrm{CH}^{(2)}_{s}\left(
\mathfrak{Q}^{(2),0,j}
\right).
$$

Since $(\mathfrak{P}^{(2)},\mathfrak{Q}^{(2)})\in\mathrm{Ker}(\mathrm{CH}^{(2)})_{s}$, we have that $(\mathfrak{P}^{(2),i+1,\bb{\mathfrak{P}^{(2)}}-1},\mathfrak{Q}^{(2),j+1,\bb{\mathfrak{Q}^{(2)}}-1})$ and $(\mathfrak{P}^{(2),0,i}, \mathfrak{Q}^{(2),0,j})$ are in $\mathrm{Ker}(\mathrm{CH}^{(2)})_{s}$. Note that, according to Definition~\ref{DDOrd}, we have that $(\mathfrak{P}^{(2),0,i},s)$ and $(\mathfrak{P}^{(2),i+1,\bb{\mathfrak{P}^{(2)}}-1},s)$ $\prec_{\mathbf{Pth}_{\boldsymbol{\mathcal{A}}^{(2)}}}$-precede $(\mathfrak{P}^{(2)},s)$. 

Therefore, by the inductive hypothesis, the second-order paths $\mathfrak{P}^{(2),0,i}$ and $\mathfrak{Q}^{(2),0,j}$, and the second-order paths $\mathfrak{P}^{(2),i+1,\bb{\mathfrak{P}^{(2)}}-1}$ and $\mathfrak{Q}^{(2),j+1,\bb{\mathfrak{Q}^{(2)}}-1}$ have, respectively, the same length, the same $([1],2)$-source and the same $([1],2)$-target. In particular $i=j$.

Therefore, we have that
\leqnomode
\allowdisplaybreaks
\begin{align*}
\bb{\mathfrak{P}^{(2)}}&=
\bb{\mathfrak{P}^{(2),0,i}}+
\bb{\mathfrak{P}^{(2),i+1,\bb{\mathfrak{P}^{(2)}}-1}}
\tag{i}
\\&=
\bb{\mathfrak{Q}^{(2),0,j}}+
\bb{\mathfrak{Q}^{(2),j+1,\bb{\mathfrak{Q}^{(2)}}-1}}
\\&=\bb{\mathfrak{Q}^{(2)}};
\\
\mathrm{sc}^{([1],2)}_{s}\left(
\mathfrak{P}^{(2)}
\right)
&=
\mathrm{sc}^{([1],2)}_{s}\left(
\mathfrak{P}^{(2),0,i}
\right)
\tag{ii}
\\&=
\mathrm{sc}^{([1],2)}_{s}\left(
\mathfrak{Q}^{(2),0,j}
\right)
\\&=
\mathrm{sc}^{([1],2)}_{s}\left(
\mathfrak{Q}^{(2)}\right)
;
\\
\mathrm{tg}^{([1],2)}_{s}\left(
\mathfrak{P}^{(2)}
\right)
&=
\mathrm{tg}^{([1],2)}_{s}\left(
\mathfrak{P}^{(2),i+1,\bb{\mathfrak{P}^{(2)}}-1}
\right)
\tag{iii}
\\&=
\mathrm{tg}^{([1],2)}_{s}\left(
\mathfrak{Q}^{(2),j+1,\bb{\mathfrak{Q}^{(2)}}-1}
\right)
\\&=
\mathrm{tg}^{([1],2)}_{s}\left(
\mathfrak{Q}^{(2)}
\right)
.
\end{align*}

The case of $\mathfrak{P}^{(2)}$ being a head-constant echelonless second-order path that is not coherent follows.

If~(2.3), then there exists a unique word $\mathbf{s}\in S^{\star}-\{\lambda\}$ and a unique operation symbol $\tau\in\Sigma^{\boldsymbol{\mathcal{A}}}_{\mathbf{s},s}$ associated to $\mathfrak{P}^{(2)}$. Let $(\mathfrak{P}^{(2)}_{j})_{j\in\bb{\mathbf{s}}}$ be the family of second-order paths in $\mathbf{Pth}_{\boldsymbol{\mathcal{A}},\mathbf{s}}$ which, in virtue of Lemma~\ref{LDPthExtract}, we can extract from $\mathfrak{P}^{(2)}$. Then, according to Definition~\ref{DDCH}, the image of the second-order Curry-Howard mapping at $\mathfrak{P}^{(2)}$ is given by
$$
\mathrm{CH}^{(2)}_{s}\left(
\mathfrak{P}^{(2)}
\right)
=
\tau^{\mathbf{T}_{\Sigma^{\boldsymbol{\mathcal{A}}^{(2)}}}(X)}
\left(\left(
\mathrm{CH}^{(2)}_{s_{j}}\left(
\mathfrak{P}^{(2)}_{j}
\right)\right)_{j\in\bb{\mathbf{s}}}
\right).
$$

Since $\mathrm{CH}^{(2)}_{s}(\mathfrak{P}^{(2)})\in\mathcal{T}(\tau,\mathrm{T}_{\Sigma^{\boldsymbol{\mathcal{A}}^{(2)}}}(X))_{1}$, which is a subset of $\mathrm{T}_{\Sigma^{\boldsymbol{\mathcal{A}}^{(2)}}}(X)^{\mathsf{HdC}\And\mathsf{C}}_{s}$, and $(\mathfrak{P}^{(2)}, \mathfrak{Q}^{(2)})\in\mathrm{Ker}(\mathrm{CH}^{(2)})_{s}$ we have, by Lemma~\ref{LDCHNEchHdC}, that $\mathfrak{Q}^{(2)}$ is a head-constant coherent echelonless second-order path associated to $\tau$, the same operation symbol as that associated to $\mathfrak{P}^{(2)}$.

Let $(\mathfrak{Q}^{(2)}_{j})_{j\in\bb{\mathbf{s}}}$ be the family of second-order paths in $\mathrm{Pth}_{\boldsymbol{\mathcal{A}}^{(2)},\mathbf{s}}$ which, by Lemma~\ref{LDPthExtract}, we can extract from $\mathfrak{Q}^{(2)}$. Then, according to Definition~\ref{DDCH}, the image of the second-order Curry-Howard mapping at $\mathfrak{Q}^{(2)}$ is given by 
$$
\mathrm{CH}^{(2)}_{s}\left(
\mathfrak{Q}^{(2)}
\right)
=
\tau^{\mathbf{T}_{\Sigma^{\boldsymbol{\mathcal{A}}^{(2)}}}(X)}
\left(\left(
\mathrm{CH}^{(2)}_{s_{j}}\left(
\mathfrak{Q}^{(2)}_{j}
\right)\right)_{j\in\bb{\mathbf{s}}}
\right).
$$

Since $(\mathfrak{P}^{(2)},\mathfrak{Q}^{(2)})\in\mathrm{Ker}(\mathrm{CH}^{(2)})_{s}$, we have, for every $j\in\bb{\mathbf{s}}$, that $(\mathfrak{P}^{(2)}_{j}, \mathfrak{Q}^{(2)}_{j})\in\mathrm{Ker}(\mathrm{CH}^{(2)})_{s_{j}}$. Note that, according to Definition~\ref{DDOrd}, we have that, for every $j\in\bb{\mathbf{s}}$, $(\mathfrak{P}^{(2)}_{j}, s_{j})$ $\prec_{\mathbf{Pth}_{\boldsymbol{\mathcal{A}}^{(2)}}}$-precedes $(\mathfrak{P}^{(2)},s)$.

Therefore, by the inductive hypothesis, for every $j\in\bb{\mathbf{s}}$, the second-order paths $\mathfrak{P}^{(2)}_{j}$ and $\mathfrak{Q}^{(2)}_{j}$ have the same length, the same $([1],2)$-source and the same $([1],2)$-target.

Therefore we have that
\leqnomode
\allowdisplaybreaks
\begin{align*}
\bb{\mathfrak{P}^{(2)}}
&=\sum_{j\in\bb{\mathbf{s}}}
\bb{\mathfrak{P}^{(2)}_{j}}
\tag{i}
\\&=\sum_{j\in\bb{\mathbf{s}}}
\bb{\mathfrak{Q}^{(2)}_{j}}
\\&=\bb{\mathfrak{Q}^{(2)}};
\\
\mathrm{sc}^{([1],2)}_{s}\left(
\mathfrak{P}^{(2)}
\right)
&=\tau^{[
\mathbf{PT}_{\boldsymbol{\mathcal{A}}}
]}\left(\left(
\mathrm{sc}^{([1],2)}_{s_{j}}\left(
\mathfrak{P}^{(2)}_{j}
\right)\right)_{j\in\bb{\mathbf{s}}}
\right)
\tag{ii}
\\&=\tau^{[
\mathbf{PT}_{\boldsymbol{\mathcal{A}}}
]}\left(\left(
\mathrm{sc}^{([1],2)}_{s_{j}}\left(
\mathfrak{Q}^{(2)}_{j}
\right)\right)_{j\in\bb{\mathbf{s}}}
\right)
\\&=
\mathrm{sc}^{([1],2)}_{s}\left(
\mathfrak{Q}^{(2)}
\right)
;
\\
\mathrm{tg}^{([1],2)}_{s}\left(
\mathfrak{P}^{(2)}
\right)
&=\tau^{[
\mathbf{PT}_{\boldsymbol{\mathcal{A}}}
]}\left(\left(
\mathrm{tg}^{([1],2)}_{s_{j}}\left(
\mathfrak{P}^{(2)}_{j}
\right)\right)_{j\in\bb{\mathbf{s}}}
\right)
\tag{iii}
\\&=\tau^{[
\mathbf{PT}_{\boldsymbol{\mathcal{A}}}
]}\left(\left(
\mathrm{tg}^{([1],2)}_{s_{j}}\left(
\mathfrak{Q}^{(2)}_{j}
\right)\right)_{j\in\bb{\mathbf{s}}}
\right)
\\&=
\mathrm{tg}^{([1],2)}_{s}\left(
\mathfrak{Q}^{(2)}
\right)
.
\end{align*}

Note that the last two items follow from the fact that $\Theta^{[1]}$ is a $\Sigma^{\boldsymbol{\mathcal{A}}}$-congruence on $\mathbf{T}_{\Sigma^{\boldsymbol{\mathcal{A}}}}(X)$ in virtue of Definition~\ref{DThetaCong}.

The case of $\mathfrak{P}^{(2)}$ being a head-constant coherent echelonless second-order path  follows.

This finishes the proof.
\end{proof}

\begin{restatable}{corollary}{CDCH}
\label{CDCH}
Let $s$ be a sort in $S$ and $(\mathfrak{P}^{(2)},\mathfrak{Q}^{(2)})\in\mathrm{Ker}(\mathrm{CH}^{(2)})_{s}$. Then the following statements hold
\begin{enumerate}
\item[(i)]  $\mathrm{sc}^{(0,2)}_{s}(\mathfrak{P}^{(2)})=\mathrm{sc}^{(0,2)}_{s}(\mathfrak{Q}^{(2)})$;
\item[(ii)] $\mathrm{tg}^{(0,2)}_{s}(\mathfrak{P}^{(2)})=\mathrm{tg}^{(0,2)}_{s}(\mathfrak{Q}^{(2)})$.
\end{enumerate}
\end{restatable}

\begin{proof}
Unravelling the mappings $\mathrm{sc}^{(0,2)}$ and $\mathrm{tg}^{(0,2)}$, introduced in Definition~\ref{DDScTgZ}, we have that 
\begin{align*}
\mathrm{sc}^{(0,2)}&=
\mathrm{sc}^{(0,[1])}
\circ
\mathrm{ip}^{([1],X)@}
\circ
\mathrm{sc}^{([1],2)};
&
\mathrm{tg}^{(0,2)}&=
\mathrm{tg}^{(0,[1])}
\circ
\mathrm{ip}^{([1],X)@}
\circ
\mathrm{tg}^{([1],2)}.
\end{align*}

Since $(\mathfrak{P}^{(2)},\mathfrak{Q}^{(2)})\in\mathrm{Ker}(\mathrm{CH}^{(2)})_{s}$, then, from Lemma~\ref{LDCH}, we have that 
\begin{align*}
\mathrm{sc}^{([1],2)}_{s}\left(
\mathfrak{P}^{(2)}
\right)&
=\mathrm{sc}^{([1],2)}_{s}\left(
\mathfrak{Q}^{(2)}
\right);
&
\mathrm{tg}^{([1],2)}_{s}\left(
\mathfrak{P}^{(2)}
\right)
&=\mathrm{tg}^{([1],2)}_{s}\left(
\mathfrak{Q}^{(2)}
\right).
\end{align*} 

The statement follows directly.
\end{proof}

Next proposition will be later used. It states that the second-order Curry-Howard mapping
restricted to the subset of $\mathrm{Pth}_{\boldsymbol{\mathcal{A}}^{(2)}}$ formed by the one-step second-order paths is an injective mapping.

\begin{restatable}{proposition}{PDCHOneStep}
\label{PDCHOneStep}
Let $s$ be a sort in $S$ and $(\mathfrak{P}^{(2)},\mathfrak{Q}^{(2)})\in\mathrm{Ker}(\mathrm{CH}^{(2)})_{s}$. If $\mathfrak{P}^{(2)}$ or $\mathfrak{Q}^{(2)}$ is a one-step second-order path, then $\mathfrak{Q}^{(2)}=\mathfrak{P}^{(2)}$.
\end{restatable}

\begin{proof}
Let $\mathfrak{P}^{(2)}$ be a one-step second-order path.
We prove the statement by Artinian induction on $(\coprod\mathrm{Pth}_{\boldsymbol{\mathcal{A}}^{(2)}},\leq_{\mathbf{Pth}_{\boldsymbol{\mathcal{A}}^{(2)}}})$.

\textsf{Base step of the Artinian induction}.

Let $(\mathfrak{P}^{(2)},s)$ be a minimal element of $(\coprod\mathrm{Pth}_{\boldsymbol{\mathcal{A}}^{(2)}}, \leq_{\mathbf{Pth}_{\boldsymbol{\mathcal{A}}}^{(2)}})$. Then, by Proposition~\ref{PDMinimal}, the second-order path $\mathfrak{P}^{(2)}$ is either~(1) a $(2,[1])$-identity second-order path on a minimal path term class, or~(2) a second-order echelon. But, since we are assuming that $\mathfrak{P}^{(2)}$ is a one-step second-order path, it follows that $\mathfrak{P}^{(2)}$ can only be a second-order echelon. This case follows from Proposition~\ref{PDCHEch}.

\textsf{Inductive step of the Artinian induction}.

Let $(\mathfrak{P}^{(2)},s)$ be a non-minimal element of $(\coprod\mathrm{Pth}_{\boldsymbol{\mathcal{A}}^{(2)}}, \leq_{\mathbf{Pth}_{\boldsymbol{\mathcal{A}}^{(2)}}})$. Let us suppose that, for every sort $t\in S$ and every second-order path $\mathfrak{P}'^{(2)}\in\mathrm{Pth}_{\boldsymbol{\mathcal{A}}^{(2)},t}$, if $(\mathfrak{P}'^{(2)},t)<_{\mathbf{Pth}_{\boldsymbol{\mathcal{A}}}^{(2)}}(\mathfrak{P}^{(2)},s)$, then the statement holds for $\mathfrak{P}'^{(2)}$, i.e., if $\mathfrak{P}'^{(2)}$ is a one-step second-order path in $\mathrm{Pth}_{\boldsymbol{\mathcal{A}}^{(2)},t}$, and $\mathfrak{Q}'^{(2)}$ is a second-order path in $\mathrm{Pth}_{\boldsymbol{\mathcal{A}}^{(2)},t}$ such that $(\mathfrak{P}'^{(2)},\mathfrak{Q}'^{(2)})\in\mathrm{Ker}(\mathrm{CH}^{(2)})_{t}$, then $\mathfrak{P}'^{(2)}=\mathfrak{Q}'^{(2)}$.

Since $(\mathfrak{P}^{(2)},s)$ is a non-minimal element in $(\coprod\mathrm{Pth}_{\boldsymbol{\mathcal{A}}^{(2)}}, \leq_{\mathbf{Pth}_{\boldsymbol{\mathcal{A}}^{(2)}}})$, it follows, by Lemma~\ref{LDOrdI}, that $\mathfrak{P}^{(2)}$ is either~(1) a $(2,[1])$-identity second-order path on a non-minimal path term class, or~(2) a second-order path of length strictly greater than one containing at least one second-order echelon or~(3) an echelonless second-order path. But, since we are assuming that $\mathfrak{P}^{(2)}$ is a one-step second-order path, it follows that $\mathfrak{P}^{(2)}$ can only be an echelonless second-order path. According to Corollary~\ref{CDCohOneStep}, $\mathfrak{P}^{(2)}$ is a coherent head-constant echelonless second-order path.

Thus, following Definition~\ref{DDHeadCt}, there exists a unique word $\mathbf{s}\in S^{\star}-\{\lambda\}$ and a unique operation symbol $\tau\in \Sigma^{\boldsymbol{\mathcal{A}}}_{\mathbf{s},s}$ associated to $\mathfrak{P}^{(2)}$. Let $(\mathfrak{P}^{(2)}_{j})_{j\in\bb{\mathbf{s}}}$ be the family of second-order paths in $\mathrm{Pth}_{\boldsymbol{\mathcal{A}}^{(2)},\mathbf{s}}$ which, in virtue of Lemma~\ref{LDPthExtract}, we can extract from $\mathfrak{P}^{(2)}$. In this case, the value of the second-order Curry-Howard mapping at $\mathfrak{P}^{(2)}$ is given by 
$$
\mathrm{CH}^{(2)}_{s}\left(
\mathfrak{P}^{(2)}
\right)=
\tau^{\mathbf{T}_{\Sigma^{\boldsymbol{\mathcal{A}}^{(2)}}}(X)}
\left(\left(\mathrm{CH}^{(2)}_{s_{j}}\left(
\mathfrak{P}^{(2)}_{j}
\right)\right)_{j\in\bb{\mathbf{s}}}
\right).
$$

Since $\bb{\mathfrak{P}^{(2)}}\geq 1$, there exists at least one $k\in\bb{\mathbf{s}}$ such that $\mathfrak{P}^{(2)}_{k}$ is a non-$(2,[1])$-identity second-order path. Since $\mathfrak{P}^{(2)}$ is a one-step second-order path, we have that $\mathfrak{P}^{(2)}_{k}$ is also a one-step second-order path. Moreover, for every $j\in\bb{\mathbf{s}}-\{k\}$, $\mathfrak{P}^{(2)}_{j}$ is a $(2,[1])$-identity second-order path. Let us note that, according to Definition~\ref{DDOrd}, $(\mathfrak{P}^{(2)}_{k},s_{k})\prec_{\mathbf{Pth}_{\boldsymbol{\mathcal{A}}^{(2)}}}(\mathfrak{P}^{(2)},s)$.

Since $(\mathfrak{P}^{(2)},\mathfrak{Q}^{(2)})\in\mathrm{Ker}(\mathrm{CH}^{(2)})_{s}$, and $\mathrm{CH}^{(2)}_{s}(\mathfrak{P}^{(2)})\in \mathcal{T}(\tau, \mathrm{T}_{\Sigma^{\boldsymbol{\mathcal{A}}^{(2)}}}(X))_{1}$, we conclude, in virtue of Lemma~\ref{LDCHNEchHdC}, that  $\mathfrak{Q}^{(2)}$ is a coherent head-constant echelonless second-order path associated to $\tau$. In fact $\mathfrak{Q}^{(2)}$ is also a one-step second-order path in virtue of Lemma~\ref{LDCH}. Let $(\mathfrak{Q}^{(2)}_{j})_{j\in\bb{\mathbf{s}}}$ be the family of second-order paths which we can extract from $\mathfrak{Q}^{(2)}$ in virtue of Lemma~\ref{LDPthExtract}. By construction, the value of the second-order Curry-Howard mapping at $\mathfrak{Q}^{(2)}$ is given by 
$$
\mathrm{CH}^{(2)}_{s}\left(
\mathfrak{Q}^{(2)}
\right)=
\tau^{\mathbf{T}_{\Sigma^{\boldsymbol{\mathcal{A}}}}(X)}
\left(\left(\mathrm{CH}^{(2)}_{s_{j}}\left(
\mathfrak{Q}^{(2)}_{j}
\right)\right)_{j\in\bb{\mathbf{s}}}
\right).
$$

Thus, for every $j\in\bb{\mathbf{s}}$, the pairs $(\mathfrak{P}^{(2)}_{j},\mathfrak{Q}^{(2)}_{j})$ are in $\mathrm{Ker}(\mathrm{CH}^{(2)})_{s_{j}}$. Then, by Corollary~\ref{CDCHUId}, we have that, for every $j\in\bb{\mathbf{s}}-\{k\}$, $\mathfrak{P}^{(2)}_{j}=\mathfrak{Q}^{(2)}_{j}$. In particular, for the index $k\in\bb{\mathbf{s}}$ we have, by the inductive hypothesis, that $\mathfrak{P}^{(2)}_{k}=\mathfrak{Q}^{(2)}_{k}$.

Since $\mathfrak{P}^{(2)}$ and $\mathfrak{Q}^{(2)}$ are one-step second-order paths, we have, by Corollary~\ref{CDUStep}, that
$$
\mathfrak{P}^{(2)}=
\tau^{\mathbf{Pth}_{\boldsymbol{\mathcal{A}}^{(2)}}}
\left(\left(\mathfrak{P}^{(2)}_{j}
\right)_{j\in\bb{\mathbf{s}}}
\right)=
\tau^{\mathbf{Pth}_{\boldsymbol{\mathcal{A}}^{(2)}}}
\left(\left(\mathfrak{Q}^{(2)}_{j}
\right)
_{j\in\bb{\mathbf{s}}}
\right)
=\mathfrak{Q}^{(2)}.
$$

This finishes the proof.
\end{proof}

We next introduce a technical lemma aiming at describing the following situation: given as inputs two  coherent head-constant echelonless second-order paths whose $1$-composition is defined but, contrarily to the inputs, the resulting $1$-composition is head-constant and echelonless but not coherent. Now consider two new second-order paths in the same class, respectively, of the original inputs with respect to the kernel of the second-order Curry-Howard mapping. Then the next lemma proves that these new second-order paths necessarily have the same behaviour as the original inputs and its $1$-composition behaves as the $1$-composition of the inputs.

\begin{restatable}{lemma}{LTech}
\label{LTech}
Let $s$ be a sort in $S$ and let $\mathfrak{P}'$ and $\mathfrak{P}$ be second-order paths in $\mathrm{Pth}_{\boldsymbol{\mathcal{A}}^{(2)},s}$ satisfying the following conditions
\begin{enumerate}
\item[(i)] $\mathrm{sc}^{([1],2)}_{s}(\mathfrak{P}'^{(2)})=\mathrm{tg}^{([1],2)}_{s}(\mathfrak{P}^{(2)})$;
\item[(ii)]  $\mathfrak{P}'^{(2)}
\circ_{s}^{1\mathbf{Pth}_{\boldsymbol{\mathcal{A}}^{(2)}}}
\mathfrak{P}^{(2)}
$ is a head-constant echelonless second-order path that is not coherent and $i=\bb{\mathfrak{P}^{(2)}}-1$ is the greatest index for which the subpath $(\mathfrak{P}'^{(2)}
\circ_{s}^{1\mathbf{Pth}_{\boldsymbol{\mathcal{A}}^{(2)}}}
\mathfrak{P}^{(2)})^{0,i}$ is coherent.
\end{enumerate}

Let $s$ be a sort in $S$ and let $\mathfrak{Q}'^{(2)}$ and $\mathfrak{Q}^{(2)}$ be second-order paths in $\mathrm{Pth}_{\boldsymbol{\mathcal{A}}^{(2)},s}$ satisfying that
$(
\mathfrak{P}'^{(2)},
\mathfrak{Q}'^{(2)}
)$,  
$(
\mathfrak{P}^{(2)},
\mathfrak{Q}^{(2)}
)
\in\mathrm{Ker}(
\mathrm{CH}^{(2)}
)_{s}.
$

Then, the following conditions holds
\begin{enumerate}
\item[(i)] $\mathrm{sc}^{([1],2)}_{s}(\mathfrak{Q}'^{(2)})=\mathrm{tg}^{([1],2)}_{s}(\mathfrak{Q}^{(2)})$;
\item[(ii)]  $\mathfrak{Q}'^{(2)}
\circ_{s}^{1\mathbf{Pth}_{\boldsymbol{\mathcal{A}}^{(2)}}}
\mathfrak{Q}^{(2)}
$ is a head-constant echelonless second-order path that is not coherent and $i=\bb{\mathfrak{Q}^{(2)}}-1$ is the greatest index for which the subpath $(\mathfrak{Q}'^{(2)}
\circ_{s}^{1\mathbf{Pth}_{\boldsymbol{\mathcal{A}}^{(2)}}}
\mathfrak{Q}^{(2)})^{0,i}$ is coherent.
\end{enumerate}
\end{restatable}
\begin{proof}

The following chain of equalities holds
\begin{align*}
\mathrm{sc}^{([1],2)}_{s}\left(
\mathfrak{Q}'^{(2)}
\right)&=
\mathrm{sc}^{([1],2)}_{s}\left(
\mathfrak{P}'^{(2)}
\right)
\tag{1}
\\&=
\mathrm{tg}^{([1],2)}_{s}\left(
\mathfrak{P}^{(2)}
\right)
\tag{2}
\\&=
\mathrm{tg}^{([1],2)}_{s}\left(
\mathfrak{Q}^{(2)}
\right).
\tag{3}
\end{align*}

The first equality follows from the fact that $(
\mathfrak{P}'^{(2)},
\mathfrak{Q}'^{(2)}
)$ is a pair in $\mathrm{Ker}(\mathrm{CH}^{(2)})_{s}$ and Lemma~\ref{LDCH}; the second equality follows by assumption; finally, the last equality follows from the fact that $(
\mathfrak{P}^{(2)},
\mathfrak{Q}^{(2)}
)$ is a pair in $\mathrm{Ker}(\mathrm{CH}^{(2)})_{s}$ and Lemma~\ref{LDCH}.

This completes the study of item (i).

Regarding item (ii), since $\mathfrak{P}'^{(2)}
\circ_{s}^{1\mathbf{Pth}_{\boldsymbol{\mathcal{A}}^{(2)}}}
\mathfrak{P}^{(2)}
$ is a head-constant echelonless second-order path that is not coherent and $i=\bb{\mathfrak{P}^{(2)}}-1$ is the greatest index for which $(\mathfrak{P}'^{(2)}
\circ_{s}^{1\mathbf{Pth}_{\boldsymbol{\mathcal{A}}^{(2)}}}
\mathfrak{P}^{(2)})^{0,i}$ is coherent, we conclude that
\begin{enumerate}
\item $\mathfrak{P}^{(2)}$ is a coherent head-constant echelonless second-order path; and
\item $\mathfrak{P}'^{(2)}$ is a head-constant echelonless second-order path.
\end{enumerate}

We first prove that, since 
$(
\mathfrak{P}'^{(2)},
\mathfrak{Q}'^{(2)}
)$ and   
$(
\mathfrak{P}^{(2)},
\mathfrak{Q}^{(2)}
)
$ are pairs in $\mathrm{Ker}(
\mathrm{CH}^{(2)}
)_{s}$, a similar situation happens for $\mathfrak{Q}'^{(2)}$ and $\mathfrak{Q}^{(2)}$, respectively.

Consider (1) the second-order path $\mathfrak{P}^{(2)}$. Then there exists a unique word $\mathbf{s}\in S^{\star}-\{\lambda\}$ and  a unique operation symbol $\tau\in \Sigma^{\boldsymbol{\mathcal{A}}}_{\mathbf{s},s}$ associated to $\mathfrak{P}^{(2)}$.
Let $(\mathfrak{P}^{(2)}_{j})_{j\in\bb{\mathbf{s}}}$ be the family of second-order paths we can extract from $\mathfrak{P}^{(2)}$ in virtue of Lemma~\ref{LDPthExtract}. Then, according to Definition~\ref{DDCH}, the value of the second-order Curry-Howard mapping at $\mathfrak{P}^{(2)}$ is given by
$$
\mathrm{CH}^{(2)}_{s}\left(
\mathfrak{P}^{(2)}
\right)
=
\tau^{\mathbf{T}_{\Sigma^{\boldsymbol{\mathcal{A}}^{(2)}}}(X)}
\left(\left(
\mathrm{CH}^{(2)}_{s_{j}}\left(
\mathfrak{P}^{(2)}_{j}
\right)\right)_{j\in\bb{\mathbf{s}}}
\right).
$$

By Lemma~\ref{LDCHNEchHdC}, we have that $\mathrm{CH}^{(2)}_{s}(\mathfrak{P}^{(2)})\in\mathcal{T}(\tau,\mathrm{T}_{\Sigma^{\boldsymbol{\mathcal{A}}^{(2)}}}(X))_{1}$, which is a subset of $\mathrm{T}_{\Sigma^{\boldsymbol{\mathcal{A}}^{(2)}}}(X)^{\mathsf{HdC}\And\mathsf{C}}_{s}$. Since $(\mathfrak{P}^{(2)},\mathfrak{Q}^{(2)})$ is in $\mathrm{Ker}(\mathrm{CH}^{(2)})_{s}$, we have, again by Lemma~\ref{LDCHNEchHdC}, that
\begin{itemize}
\item[(1)] $\mathfrak{Q}^{(2)}$ is a  coherent head-constant echelonless second-order path  associated to the operation symbol $\tau\in\Sigma^{\boldsymbol{\mathcal{A}}}_{\mathbf{s},s}$.
\end{itemize}

Let $(\mathfrak{Q}^{(2)}_{j})_{j\in\bb{\mathbf{s}}}$ be the family of second-order paths we can extract from $\mathfrak{Q}^{(2)}$ in virtue of Lemma~\ref{LDPthExtract}. Then, following Definition~\ref{DDCH}, the value of the second-order Curry-Howard mapping at $\mathfrak{Q}^{(2)}$ is given by
$$
\mathrm{CH}^{(2)}_{s}\left(
\mathfrak{Q}^{(2)}
\right)
=
\tau^{\mathbf{T}_{\Sigma^{\boldsymbol{\mathcal{A}}^{(2)}}}(X)}
\left(\left(
\mathrm{CH}^{(2)}_{s_{j}}\left(
\mathfrak{Q}^{(2)}_{j}
\right)\right)_{j\in\bb{\mathbf{s}}}
\right).
$$

Since $(\mathfrak{P}^{(2)},\mathfrak{Q}^{(2)})$ is a pair in $\mathrm{Ker}(\mathrm{CH}^{(2)})_{s}$, we have, for every $j\in\bb{\mathbf{s}}$, that
$$
\left(\mathfrak{P}^{(2)}_{j},
\mathfrak{Q}^{(2)}_{j}
\right)
\in\mathrm{Ker}\left(\mathrm{CH}^{(2)}\right)_{s_{j}}.
$$

This completes the study of item (1).

Now, consider (2) the second-order path $\mathfrak{P}'^{(2)}$. Let us note that $\mathfrak{P}'^{(2)}$ has to be a head-constant echelonless second-order path associated to the operation symbol $\tau$. Note that the operation symbol $\tau$ is the same as in case (1), since $\mathfrak{P}'^{(2)}\circ_{s}^{1\mathbf{Pth}_{\boldsymbol{\mathcal{A}}^{(2)}}}\mathfrak{P}^{(2)}$ is head-constant by item (ii) in the hypothesis. 

Now, it could be the case that $\mathfrak{P}'^{(2)}$ is either coherent or not. We will assume that $\mathfrak{P}'^{(2)}$ is not coherent. The coherent case will be handled similarly. Since $\mathfrak{P}'^{(2)}$ is a head-constant echelonless second-order path of length at least one that is not coherent, let $k\in\bb{\mathfrak{P}'^{(2)}}$ be the greatest index for which $\mathfrak{P}'^{(2),0,k}$ is coherent.

According to Definition~\ref{DDCH}, we have that the value of the second-order Curry-Howard mapping at $\mathfrak{P}'^{(2)}$ is given by
$$
\mathrm{CH}^{(2)}_{s}\left(
\mathfrak{P}'^{(2)}
\right)
=
\mathrm{CH}^{(2)}_{s}\left(
\mathfrak{P}'^{(2), k+1,\bb{\mathfrak{P}'^{(2)}}-1}
\right)
\circ_{s}^{1\mathbf{T}_{\Sigma^{\boldsymbol{\mathcal{A}}^{(2)}}}(X)}
\mathrm{CH}^{(2)}_{s}\left(
\mathfrak{P}'^{(2),0,k}
\right).
$$

By Lemma~\ref{LDCHNEchHd}, we have that $\mathrm{CH}^{(2)}_{s}(\mathfrak{P}'^{(2)})\in\mathcal{T}(\tau,\mathrm{T}_{\Sigma^{\boldsymbol{\mathcal{A}}^{(2)}}}(X))^{+}$, which is a subset of $\mathrm{T}_{\Sigma^{\boldsymbol{\mathcal{A}}^{(2)}}}(X)^{\mathsf{HdC}\And\neg\mathsf{C}}_{s}$. Since $(\mathfrak{P}'^{(2)},\mathfrak{Q}'^{(2)})$ is in $\mathrm{Ker}(\mathrm{CH}^{(2)})_{s}$, we have, again by Lemma~\ref{LDCHNEchHdC}, that
\begin{itemize}
\item[(2)] $\mathfrak{Q}'^{(2)}$ is a head-constant echelonless second-order path that is not coherent  associated to the operation symbol $\tau\in\Sigma^{\boldsymbol{\mathcal{A}}}_{\mathbf{s},s}$.
\end{itemize}

Let $l\in\bb{\mathfrak{Q}'^{(2)}}$ be the greatest index for which $\mathfrak{Q}'^{(2),0,l}$ is coherent. According to Definition~\ref{DDCH}, we have that the value of the second-order Curry-Howard mapping at $\mathfrak{Q}'^{(2)}$ is given by
$$
\mathrm{CH}^{(2)}_{s}\left(
\mathfrak{Q}'^{(2)}
\right)
=
\mathrm{CH}^{(2)}_{s}\left(
\mathfrak{Q}'^{(2),l+1,\bb{\mathfrak{Q}'^{(2)}}-1}
\right)
\circ_{s}^{1\mathbf{T}_{\Sigma^{\boldsymbol{\mathcal{A}}^{(2)}}}(X)}
\mathrm{CH}^{(2)}_{s}\left(
\mathfrak{Q}'^{(2),0,l}
\right).
$$

Since $(\mathfrak{P}'^{(2)},\mathfrak{Q}'^{(2)})\in\mathrm{Ker}(\mathrm{CH}^{(2)})_{s}$, we conclude that $(\mathfrak{P}'^{(2),k+1,\bb{\mathfrak{P}'^{(2)}}-1}, \mathfrak{Q}'^{(2),l+1,\bb{\mathfrak{Q}'^{(2)}}-1})$ and $(\mathfrak{P}'^{(2),0,k},\mathfrak{Q}'^{(2),0,l})$ are in $\mathrm{Ker}(\mathrm{CH}^{(2)})_{s}$. 

Now, back to $\mathfrak{P}'^{(2),0,k}$. This is a coherent head-constant echelonless second-order path, let $(\mathfrak{P}'^{(2)}_{j})_{j\in\bb{\mathbf{s}}}$ be the family of second-order paths we can extract from $\mathfrak{P}'^{(2),0,k}$ in virtue of Lemma~\ref{LDPthExtract}. Then, according to Definition~\ref{DDCH}, we have that the value of the second-order Curry-Howard mapping at $\mathfrak{P}'^{(2),0,k}$ is given by
$$
\mathrm{CH}^{(2)}_{s}\left(
\mathfrak{P}'^{(2),0,k}
\right)
=
\tau^{\mathbf{T}_{\Sigma^{\boldsymbol{\mathcal{A}}^{(2)}}}(X)}
\left(\left(
\mathrm{CH}^{(2)}_{s_{j}}\left(
\mathfrak{P}'^{(2)}_{j}
\right)\right)_{j\in\bb{\mathbf{s}}}
\right).
$$

Note that $\mathrm{CH}^{(2)}_{s}(\mathfrak{P}'^{(2),0,k})$ is a term in $\mathcal{T}(\tau,\mathrm{T}_{\Sigma^{\boldsymbol{\mathcal{A}}^{(2)}}}(X))_{1}$, which is a subset of $\mathrm{T}_{\Sigma^{\boldsymbol{\mathcal{A}}^{(2)}}}(X)^{\mathsf{HdC}\And\mathsf{C}}_{s}$. Since $(\mathfrak{P}'^{(2),0,k},\mathfrak{Q}'^{(2),0,l})$ is in $\mathrm{Ker}(\mathrm{CH}^{(2)})_{s}$, following Lemma~\ref{LDCHNEchHdC}, we conclude that $\mathfrak{Q}'^{(2),0,l}$ is a coherent head-constant echelonless second-order path associated to the operation symbol $\tau$.

Let $(\mathfrak{Q}'^{(2)}_{j})_{j\in\bb{\mathbf{s}}}$ be the family of second-order paths we can extract from $\mathfrak{Q}'^{(2),0,l}$ in virtue of Lemma~\ref{LDPthExtract}. Then, following Definition~\ref{DDCH}, the value of the second-order Curry-Howard mapping at $\mathfrak{Q}'^{(2)}$ is given by
$$
\mathrm{CH}^{(2)}_{s}\left(
\mathfrak{Q}'^{(2),0,l}
\right)
=
\tau^{\mathbf{T}_{\Sigma^{\boldsymbol{\mathcal{A}}^{(2)}}}(X)}
\left(\left(
\mathrm{CH}^{(2)}_{s_{j}}\left(
\mathfrak{Q}'^{(2)}_{j}
\right)\right)_{j\in\bb{\mathbf{s}}}
\right).
$$

Since $(\mathfrak{P}'^{(2),0,k},\mathfrak{Q}'^{(2),0,l})$ is a pair in $\mathrm{Ker}(\mathrm{CH}^{(2)})_{s}$, we have, for every $j\in\bb{\mathbf{s}}$, that
$$
\left(\mathfrak{P}'^{(2)}_{j},
\mathfrak{Q}'^{(2)}_{j}
\right)
\in\mathrm{Ker}\left(\mathrm{CH}^{(2)}\right)_{s_{j}}.
$$

This completes the study of item (2).

We are now in position to study the last item.

Let us begin by studying the $([1],2)$-target of $\mathfrak{P}^{(2)}$. Let us assume that $\mathfrak{P}^{(2)}$ is a second-order $\mathbf{c}$-path, for some word $\mathbf{c}\in S^{\star}-\{\lambda\}$, i.e.,
$$
\mathfrak{P}^{(2)}
=
\left(
([P_{i}]_{s})_{i\in\bb{\mathbf{c}}+1},
(\mathfrak{p}^{(2)}_{i})_{i\in\bb{\mathbf{c}}},
(T^{(1)}_{i})_{i\in\bb{\mathbf{c}}}
\right).
$$

Let us consider $\mathfrak{p}^{(2)}_{\bb{\mathbf{c}}-1}$, the last second-order rewrite rule appearing in $\mathfrak{P}^{(2)}$, and let us assume that it has the form 
$([M_{\bb{\mathbf{c}}-1}]_{c_{\bb{\mathbf{c}}-1}}
,[N_{\bb{\mathbf{c}}-1}]_{c_{\bb{\mathbf{c}}-1}})$ 
for some suitable pair of path terms classes $[M_{\bb{\mathbf{c}}-1}]_{c_{\bb{\mathbf{c}}-1}}$ and $[N_{\bb{\mathbf{c}}-1}]_{c_{\bb{\mathbf{c}}-1}}$ in $[\mathrm{PT}_{\boldsymbol{\mathcal{A}}}]_{c_{\bb{\mathbf{c}}-1}}$. Let $T^{(1)}_{\bb{\mathbf{c}}-1}$ be the last first-order translation occurring in $\mathfrak{P}^{(2)}$.

Since $\mathfrak{P}^{(2)}$ is a second-order path, we have that 
\begin{itemize}
\item[(i)]
$
T^{(1)}_{\bb{\mathbf{c}}-1}(
M_{\bb{\mathbf{c}}-1}
)
\in [P_{\bb{\mathbf{c}}-1}]_{s};
$
\item[(ii)]
$
T^{(1)}_{\bb{\mathbf{c}}-1}(
N_{\bb{\mathbf{c}}-1}
)
\in [P_{\bb{\mathbf{c}}}]_{s}.
$
\end{itemize}

Let us recall that $T^{(1)}_{\bb{\mathbf{c}}-1}$ is a first-order translation of type $\tau$. Thus there will exists an index $j_{1}\in \bb{\mathbf{s}}$,  a family of paths $(\mathfrak{P}_{k})_{k\in j_{1}}\in \prod_{k\in j_{1}}\mathrm{Pth}_{\boldsymbol{\mathcal{A}},s_{k}}$, a family of paths $(\mathfrak{P}_{l})_{l\in \bb{\mathbf{s}}-(j_{1}+1)}\in \prod_{l\in \bb{\mathbf{s}}-(j_{1}+1)}\mathrm{Pth}_{\boldsymbol{\mathcal{A}},s_{l}}$ and a first-order translation $T^{(1)'}_{\bb{\mathbf{c}}-1}\in\mathrm{Tl}_{c_{\bb{\mathbf{c}}-1}}(\mathbf{T}_{\Sigma^{\boldsymbol{\mathcal{A}}}}(X))_{s_{j_{1}}}$ such that 
\begin{multline*}
T^{(1)}_{\bb{\mathbf{c}}-1}=\tau^{\mathbf{T}_{\Sigma^{\boldsymbol{\mathcal{A}}}}(X)}\left(
\mathrm{CH}^{(1)}_{s_{0}}\left(\mathfrak{P}_{0}\right),
\cdots,
\mathrm{CH}^{(1)}_{s_{j_{1}-1}}\left(\mathfrak{P}_{j_{1}-1}\right),
\right.
\\
\left.
T^{(1)'}_{\bb{\mathbf{c}}-1},
\mathrm{CH}^{(1)}_{s_{j_{1}+1}}\left(\mathfrak{P}_{j_{1}+1}\right),
\cdots,
\mathrm{CH}^{(1)}_{s_{\bb{\mathbf{s}}-1}}\left(\mathfrak{P}_{\bb{\mathbf{s}}-1}\right)
\right).
\end{multline*}

Moreover, the following chain of equalities holds
\begin{flushleft}
$\left[
P_{\bb{\mathbf{c}}}
\right]_{s}$
\allowdisplaybreaks
\begin{align*}
\,
&=
\left[
T^{(1)}_{\bb{\mathbf{c}}-1}\left(
N_{\bb{\mathbf{c}}-1}
\right)
\right]_{s}
\tag{1}
\\&=
\biggr[
\tau^{\mathbf{PT}_{\boldsymbol{\mathcal{A}}}}
\biggr(
\mathrm{CH}^{(1)}_{s_{0}}\left(\mathfrak{P}_{0}\right),
\cdots,
\mathrm{CH}^{(1)}_{s_{j_{1}-1}}\left(\mathfrak{P}_{j_{1}-1}\right),
\\&\qquad\qquad\qquad
T^{(1)'}_{\bb{\mathbf{c}}-1}\left(
N_{\bb{\mathbf{c}}-1}
\right),
\mathrm{CH}^{(1)}_{s_{j_{1}+1}}\left(\mathfrak{P}_{j_{1}+1}\right),
\cdots,
\mathrm{CH}^{(1)}_{s_{\bb{\mathbf{s}}-1}}\left(\mathfrak{P}_{\bb{\mathbf{s}}-1}\right)
\biggr)
\biggr]_{s}
\tag{2}
\\&=
\tau^{
[\mathbf{PT}_{\boldsymbol{\mathcal{A}}}]
}
\Biggr(
\left[
\mathrm{CH}^{(1)}_{s_{0}}\left(\mathfrak{P}_{0}\right)
\right]_{s_{0}},
\cdots,
\left[
\mathrm{CH}^{(1)}_{s_{j_{1}-1}}\left(\mathfrak{P}_{j_{1}-1}\right)
\right]_{s_{j_{1}-1}},
\\&\qquad\qquad
\left[
T^{(1)'}_{\bb{\mathbf{c}}-1}\left(
N_{\bb{\mathbf{c}}-1}
\right)
\right]_{s_{j_{1}}},
\\&\qquad\qquad\qquad\qquad
\left[
\mathrm{CH}^{(1)}_{s_{j_{1}+1}}\left(\mathfrak{P}_{j_{1}+1}\right)
\right]_{s_{j_{1}+1}},
\cdots,
\left[
\mathrm{CH}^{(1)}_{s_{\bb{\mathbf{s}}-1}}\left(\mathfrak{P}_{\bb{\mathbf{s}}-1}\right)
\right]_{s_{\bb{\mathbf{s}}-1}}
\Biggr).
\tag{3}
\end{align*}
\end{flushleft}

The first equality follows from item (ii) above; the second equality unravels the description of the last first-order translation occurring in $\mathfrak{P}^{(2)}$; finally, the last equality recovers the interpretation of the operation symbol $\tau$ in the many-sorted partial $\Sigma^{\boldsymbol{\mathcal{A}}}$-algebra $[\mathbf{PT}_{\boldsymbol{\mathcal{A}}}]$, as described in Proposition~\ref{PPTQCatAlg}. 

By Lemma~\ref{LDPthExtract}, we have the following family of equalities
\allowdisplaybreaks
\begin{align*}
\mathrm{tg}^{([1],2)}_{s_{0}}\left(
\mathfrak{P}^{(2)}_{0}
\right)
&=
\left[
\mathrm{CH}^{(1)}_{s_{0}}\left(\mathfrak{P}_{0}\right)
\right]_{s_{0}};
\\
\vdots\qquad&=\qquad\vdots
\\
\mathrm{tg}^{([1],2)}_{s_{j_{1}}}
\left(
\mathfrak{P}^{(2)}_{j_{1}}
\right)
&=
\left[
T^{(1)'}_{\bb{\mathbf{c}}-1}\left(
N_{\bb{\mathbf{c}}-1}
\right)
\right]_{s_{j_{1}}};
\\
\vdots\qquad&=\qquad\vdots
\\
\mathrm{tg}^{([1],2)}_{s_{\bb{\mathbf{s}}-1}}
\left(
\mathfrak{P}^{(2)}_{s_{\bb{\mathbf{s}}-1}}
\right)
&=
\left[
\mathrm{CH}^{(1)}_{s_{\bb{\mathbf{s}}-1}}\left(\mathfrak{P}_{\bb{\mathbf{s}}-1}\right)
\right]_{s_{\bb{\mathbf{s}}-1}}.
\end{align*}

Now let us turn ourselves to the study of the $([1],2)$-target of $\mathfrak{Q}^{(2)}$. Let us assume that $\mathfrak{Q}^{(2)}$ is a second-order $\mathbf{d}$-path, for some word $\mathbf{d}\in S^{\star}$,
i.e.,
$$
\mathfrak{Q}^{(2)}
=
\left(
([Q_{i}]_{s})_{i\in\bb{\mathbf{d}}+1},
(\mathfrak{q}^{(2)}_{i})_{i\in\bb{\mathbf{d}}},
(U^{(1)}_{i})_{i\in\bb{\mathbf{d}}}
\right).
$$

Let us consider $\mathfrak{q}^{(2)}_{\bb{\mathbf{d}}-1}$, the last second-order rewrite rule appearing in $\mathfrak{Q}^{(2)}$, and let us assume that it has the form $([K_{\bb{\mathbf{d}}-1}]_{d_{\bb{\mathbf{d}}-1}},[L_{\bb{\mathbf{d}}-1}]_{d_{\bb{\mathbf{d}}-1}})$ for some suitable pair of path terms classes $[K_{\bb{\mathbf{d}}-1}]_{d_{\bb{\mathbf{d}}-1}}$ and $[L_{\bb{\mathbf{d}}-1}]_{d_{\bb{\mathbf{d}}-1}}$ in $[\mathrm{PT}_{\boldsymbol{\mathcal{A}}}]_{d_{\bb{\mathbf{d}}-1}}$. Let $U^{(1)}_{\bb{\mathbf{d}}-1}$ be the last first-order translation occurring in $\mathfrak{Q}^{(2)}$.

Since $\mathfrak{Q}^{(2)}$ is a second-order path, we have that 
\begin{itemize}
\item[(i)]
$
U^{(1)}_{\bb{\mathbf{d}}-1}(
K_{\bb{\mathbf{d}}-1}
)
\in [Q_{\bb{\mathbf{d}}-1}]_{s};
$
\item[(ii)]
$
U^{(1)}_{\bb{\mathbf{d}}-1}(
L_{\bb{\mathbf{d}}-1}
)
\in [Q_{\bb{\mathbf{d}}}]_{s}.
$
\end{itemize}

Let us recall that $U^{(1)}_{\bb{\mathbf{c}}-1}$ is a first-order translation of type $\tau$. Thus there will exists an index $j_{2}\in \bb{\mathbf{s}}$,  a family of paths $(\mathfrak{Q}_{k})_{k\in j_{2}}\in \prod_{k\in j_{2}}\mathrm{Pth}_{\boldsymbol{\mathcal{A}},s_{k}}$, a family of paths $(\mathfrak{Q}_{l})_{l\in \bb{\mathbf{s}}-(j_{2}+1)}\in \prod_{l\in \bb{\mathbf{s}}-(j_{2}+1)}\mathrm{Pth}_{\boldsymbol{\mathcal{A}},s_{l}}$ and a first-order translation $U'_{\bb{\mathbf{d}}-1}\in\mathrm{Tl}_{d_{\bb{\mathbf{d}}-1}}(\mathbf{T}_{\Sigma^{\boldsymbol{\mathcal{A}}}}(X))_{s_{j_{2}}}$ such that 
\begin{multline*}
U^{(1)}_{\bb{\mathbf{d}}-1}=\tau^{\mathbf{T}_{\Sigma^{\boldsymbol{\mathcal{A}}}}(X)}\left(
\mathrm{CH}^{(1)}_{s_{0}}\left(\mathfrak{Q}_{0}\right),
\cdots,
\mathrm{CH}^{(1)}_{s_{j_{2}-1}}\left(\mathfrak{Q}_{j_{2}-1}\right),
\right.
\\
\left.
U^{(1)'}_{\bb{\mathbf{d}}-1},
\mathrm{CH}^{(1)}_{s_{j_{2}+1}}\left(\mathfrak{Q}_{j_{2}+1}\right),
\cdots,
\mathrm{CH}^{(1)}_{s_{\bb{\mathbf{s}}-1}}\left(\mathfrak{Q}_{\bb{\mathbf{s}}-1}\right)
\right).
\end{multline*}

Moreover, the following chain of equalities holds
\begin{flushleft}
$\left[
Q_{\bb{\mathbf{d}}}
\right]_{s}$
\allowdisplaybreaks
\begin{align*}
\,
&=
\left[
U^{(1)}_{\bb{\mathbf{d}}-1}\left(
L_{\bb{\mathbf{d}}-1}
\right)
\right]_{s}
\tag{1}
\\&=
\biggr[
\tau^{\mathbf{PT}_{\boldsymbol{\mathcal{A}}}}
\biggr(
\mathrm{CH}^{(1)}_{s_{0}}\left(\mathfrak{Q}_{0}\right),
\cdots,
\mathrm{CH}^{(1)}_{s_{j_{2}-1}}\left(\mathfrak{Q}_{j_{2}-1}\right),
\\&\qquad\qquad\qquad
U^{(1)'}_{\bb{\mathbf{d}}-1}\left(
L_{\bb{\mathbf{d}}-1}
\right),
\mathrm{CH}^{(1)}_{s_{j_{2}+1}}\left(\mathfrak{Q}_{j_{2}+1}\right),
\cdots,
\mathrm{CH}^{(1)}_{s_{\bb{\mathbf{s}}-1}}\left(\mathfrak{Q}_{\bb{\mathbf{s}}-1}\right)
\biggr)
\biggr]_{s}
\tag{2}
\\&=
\tau^{
[\mathbf{PT}_{\boldsymbol{\mathcal{A}}}]
}
\Biggr(
\left[
\mathrm{CH}^{(1)}_{s_{0}}\left(\mathfrak{Q}_{0}\right)
\right]_{s_{0}},
\cdots,
\left[
\mathrm{CH}^{(1)}_{s_{j_{2}-1}}\left(\mathfrak{Q}_{j_{2}-1}\right)
\right]_{s_{j_{2}-1}},
\\&\qquad\qquad
\left[
U^{(1)'}_{\bb{\mathbf{d}}-1}\left(
L_{\bb{\mathbf{d}}-1}
\right)
\right]_{s_{j_{2}}},
\\&\qquad\qquad\qquad\qquad
\left[
\mathrm{CH}^{(1)}_{s_{j_{2}+1}}\left(\mathfrak{Q}_{j_{2}+1}\right)
\right]_{s_{j_{2}+1}},
\cdots,
\left[
\mathrm{CH}^{(1)}_{s_{\bb{\mathbf{s}}-1}}\left(\mathfrak{Q}_{\bb{\mathbf{s}}-1}\right)
\right]_{s_{\bb{\mathbf{s}}-1}}
\Biggr).
\tag{3}
\end{align*}
\end{flushleft}

The first equality follows from item (ii) above; the second equality unravels the description of the last first-order translation occurring in $\mathfrak{Q}^{(2)}$; finally, the last equality recovers the interpretation of the operation symbol $\tau$ in the many-sorted partial $\Sigma^{\boldsymbol{\mathcal{A}}}$-algebra $[\mathbf{PT}_{\boldsymbol{\mathcal{A}}}]$, as described in Proposition~\ref{PPTQCatAlg}. 

By Lemma~\ref{LDPthExtract}, we have the following family of equalities
\allowdisplaybreaks
\begin{align*}
\mathrm{tg}^{([1],2)}_{s_{0}}\left(
\mathfrak{Q}^{(2)}_{0}
\right)
&=
\left[
\mathrm{CH}^{(1)}_{s_{0}}\left(\mathfrak{Q}_{0}\right)
\right]_{s_{0}};
\\
\vdots\qquad&=\qquad\vdots
\\
\mathrm{tg}^{([1],2)}_{s_{j_{2}}}
\left(
\mathfrak{Q}^{(2)}_{j_{2}}
\right)
&=
\left[
U^{(1)'}_{\bb{\mathbf{d}}-1}\left(
L_{\bb{\mathbf{d}}-1}
\right)
\right]_{s_{j_{2}}};
\\
\vdots\qquad&=\qquad\vdots
\\
\mathrm{tg}^{([1],2)}_{s_{\bb{\mathbf{s}}-1}}
\left(
\mathfrak{Q}^{(2)}_{s_{\bb{\mathbf{s}}-1}}
\right)
&=
\left[
\mathrm{CH}^{(1)}_{s_{\bb{\mathbf{s}}-1}}\left(\mathfrak{Q}_{\bb{\mathbf{s}}-1}\right)
\right]_{s_{\bb{\mathbf{s}}-1}}.
\end{align*}

Since we have already proved that, for every $j\in\bb{\mathbf{s}}$, the pair $(\mathfrak{P}^{(2)}_{j}, \mathfrak{Q}^{(2)}_{j})$ is in $\mathrm{Ker}(\mathrm{CH}^{(2)})_{s_{j}}$, then we have, by Lemma~\ref{LDCH}, that, for every $j\in\bb{\mathbf{s}}$, it is the case that 
$$
\mathrm{tg}^{([1],2)}_{s_{j}}
\left(
\mathfrak{P}^{(2)}_{j}
\right)
=
\mathrm{tg}^{([1],2)}_{s_{j}}
\left(
\mathfrak{Q}^{(2)}_{j}
\right).
$$

Therefore, we have the following family of equalities between path term classes
\allowdisplaybreaks
\begin{align*}
\left[
\mathrm{CH}^{(1)}_{s_{0}}\left(\mathfrak{P}_{0}\right)
\right]_{s_{0}}
&=
\left[
\mathrm{CH}^{(1)}_{s_{0}}\left(\mathfrak{Q}_{0}\right)
\right]_{s_{0}};
\\
&\vdots
\\
\left[
T^{(1)'}_{\bb{\mathbf{c}}-1}\left(
N_{\bb{\mathbf{c}}-1}
\right)
\right]_{s_{j_{1}}}
&=
\left[
\mathrm{CH}^{(1)}_{s_{j_{1}}}\left(\mathfrak{Q}_{j_{1}}\right)
\right]_{s_{j_{1}}};
\\
&\vdots
\\
\left[
\mathrm{CH}^{(1)}_{s_{j_{2}}}\left(\mathfrak{P}_{j_{2}}\right)
\right]_{s_{j_{2}}}
&=
\left[
U^{(1)'}_{\bb{\mathbf{d}}-1}\left(
L_{\bb{\mathbf{d}}-1}
\right)
\right]_{s_{j_{2}}};
\\
&\vdots
\\
\left[
\mathrm{CH}^{(1)}_{s_{\bb{\mathbf{s}}-1}}\left(\mathfrak{P}_{\bb{\mathbf{s}}-1}\right)
\right]_{s_{\bb{\mathbf{s}}-1}}
&=
\left[
\mathrm{CH}^{(1)}_{s_{\bb{\mathbf{s}}-1}}\left(\mathfrak{Q}_{\bb{\mathbf{s}}-1}\right)
\right]_{s_{\bb{\mathbf{s}}-1}}.
\end{align*}

Now let us turn ourselves to the study of the $([1],2)$-source of $\mathfrak{P}'^{(2)}$. Let us assume that $\mathfrak{P}'^{(2)}$ is a second-order $\mathbf{e}$-path, for some word $\mathbf{e}\in S^{\star}-\{\lambda\}$, i.e.,
$$
\mathfrak{P}'^{(2)}
=
\left(
([R_{i}]_{s})_{i\in\bb{\mathbf{e}}+1},
(\mathfrak{r}^{(2)}_{i})_{i\in\bb{\mathbf{e}}},
(V^{(1)}_{i})_{i\in\bb{\mathbf{e}}}
\right).
$$

Let us consider $\mathfrak{r}^{(2)}_{0}$, the first second-order rewrite rule appearing in $\mathfrak{P}'^{(2)}$, and let us assume that it has the form 
$([M'_{0}]_{e_{0}}
,[N'_{0}]_{e_{0}})$ 
for some suitable pair of path terms classes $[M'_{0}]_{e_{0}}$ and $[N'_{0}]_{e_{0}}$ in $[\mathrm{PT}_{\boldsymbol{\mathcal{A}}}]_{e_{0}}$. Let $V^{(1)}_{0}$ be the initial first-order translation occurring in $\mathfrak{P}'^{(2)}$.

Since $\mathfrak{P}'^{(2)}$ is a second-order path, we have that 
\begin{itemize}
\item[(i)]
$
V^{(1)}_{0}(
M'_{0}
)
\in [R_{0}]_{s};
$
\item[(ii)]
$
V^{(1)}_{0}(
N'_{0}
)
\in [R_{1}]_{s}.
$
\end{itemize}

Let us recall that $V_{0}$ is a first-order translation of type $\tau$. Thus there will exists an index $j_{3}\in \bb{\mathbf{s}}$,  a family of paths $(\mathfrak{P}'_{k})_{k\in j_{3}}\in \prod_{k\in j_{3}}\mathrm{Pth}_{\boldsymbol{\mathcal{A}},s_{k}}$, a family of paths $(\mathfrak{P}'_{l})_{l\in \bb{\mathbf{s}}-(j_{3}+1)}\in \prod_{l\in \bb{\mathbf{s}}-(j_{3}+1)}\mathrm{Pth}_{\boldsymbol{\mathcal{A}},s_{l}}$ and a first-order translation $V^{(1)'}_{0}\in\mathrm{Tl}_{e_{0}}(\mathbf{T}_{\Sigma^{\boldsymbol{\mathcal{A}}}}(X))_{s_{j_{3}}}$ such that 
\begin{multline*}
V^{(1)}_{0}=\tau^{\mathbf{T}_{\Sigma^{\boldsymbol{\mathcal{A}}}}(X)}\left(
\mathrm{CH}^{(1)}_{s_{0}}\left(\mathfrak{P}'_{0}\right),
\cdots,
\mathrm{CH}^{(1)}_{s_{j_{3}-1}}\left(\mathfrak{P}'_{j_{3}-1}\right),
\right.
\\
\left.
V^{(1)'}_{0},
\mathrm{CH}^{(1)}_{s_{j_{3}+1}}\left(\mathfrak{P}'_{j_{3}+1}\right),
\cdots,
\mathrm{CH}^{(1)}_{s_{\bb{\mathbf{s}}-1}}\left(\mathfrak{P}'_{\bb{\mathbf{s}}-1}\right)
\right).
\end{multline*}

Moreover, the following chain of equalities holds
\begin{flushleft}
$\left[
R_{0}
\right]_{s}$
\allowdisplaybreaks
\begin{align*}
\,
&=
\left[
V^{(1)}_{0}\left(
M'_{0}
\right)
\right]_{s}
\tag{1}
\\&=
\biggr[
\tau^{\mathbf{PT}_{\boldsymbol{\mathcal{A}}}}
\biggr(
\mathrm{CH}^{(1)}_{s_{0}}\left(\mathfrak{P}'_{0}\right),
\cdots,
\mathrm{CH}^{(1)}_{s_{j_{3}-1}}\left(\mathfrak{P}'_{j_{3}-1}\right),
\\&\qquad\qquad\qquad
V^{(1)'}_{0}\left(
M'_{0}
\right),
\mathrm{CH}^{(1)}_{s_{j_{3}+1}}\left(\mathfrak{P}'_{j_{3}+1}\right),
\cdots,
\mathrm{CH}^{(1)}_{s_{\bb{\mathbf{s}}-1}}\left(\mathfrak{P}'_{\bb{\mathbf{s}}-1}\right)
\biggr)
\biggr]_{s}
\tag{2}
\\&=
\tau^{
[\mathbf{PT}_{\boldsymbol{\mathcal{A}}}]
}
\Biggr(
\left[
\mathrm{CH}^{(1)}_{s_{0}}\left(\mathfrak{P}'_{0}\right)
\right]_{s_{0}},
\cdots,
\left[
\mathrm{CH}^{(1)}_{s_{j_{3}-1}}\left(\mathfrak{P}'_{j_{3}-1}\right)
\right]_{s_{j_{3}-1}},
\\&\qquad\qquad
\left[
V^{(1)'}_{0}\left(
M'_{0}
\right)
\right]_{s_{j_{3}}},
\\&\qquad\qquad\qquad\quad
\left[
\mathrm{CH}^{(1)}_{s_{j_{3}+1}}\left(\mathfrak{P}'_{j_{3}+1}\right)
\right]_{s_{j_{3}+1}},
\cdots,
\left[
\mathrm{CH}^{(1)}_{s_{\bb{\mathbf{s}}-1}}\left(\mathfrak{P}'_{\bb{\mathbf{s}}-1}\right)
\right]_{s_{\bb{\mathbf{s}}-1}}
\Biggr).
\tag{3}
\end{align*}
\end{flushleft}

The first equality follows from item (ii) above; the second equality unravels the description of the initial first-order translation occurring in $\mathfrak{P}'^{(2)}$; finally, the last equality recovers the interpretation of the operation symbol $\tau$ in the many-sorted partial $\Sigma^{\boldsymbol{\mathcal{A}}}$-algebra $[\mathbf{PT}_{\boldsymbol{\mathcal{A}}}]$, as described in Proposition~\ref{PPTQCatAlg}. 

By Lemma~\ref{LDPthExtract}, we have the following family of equalities
\allowdisplaybreaks
\begin{align*}
\mathrm{sc}^{([1],2)}_{s_{0}}\left(
\mathfrak{P}'^{(2)}_{0}
\right)
&=
\left[
\mathrm{CH}^{(1)}_{s_{0}}\left(\mathfrak{P}'_{0}\right)
\right]_{s_{0}};
\\
\vdots\qquad&=\qquad\vdots
\\
\mathrm{sc}^{([1],2)}_{s_{j_{3}}}
\left(
\mathfrak{P}'^{(2)}_{j_{3}}
\right)
&=
\left[
V^{(1)'}_{0}\left(
M'_{0}
\right)
\right]_{s_{j_{3}}};
\\
\vdots\qquad&=\qquad\vdots
\\
\mathrm{sc}^{([1],2)}_{s_{\bb{\mathbf{s}}-1}}
\left(
\mathfrak{P}'^{(2)}_{s_{\bb{\mathbf{s}}-1}}
\right)
&=
\left[
\mathrm{CH}^{(1)}_{s_{\bb{\mathbf{s}}-1}}\left(\mathfrak{P}'_{\bb{\mathbf{s}}-1}\right)
\right]_{s_{\bb{\mathbf{s}}-1}}.
\end{align*}

We finally study of the $([1],2)$-source of $\mathfrak{Q}'^{(2)}$. Let us assume that $\mathfrak{Q}'^{(2)}$ is a second-order $\mathbf{f}$-path, for some word $\mathbf{f}\in S^{\star}-\{\lambda\}$, i.e.,
$$
\mathfrak{Q}'^{(2)}
=
\left(
([S_{i}]_{s})_{i\in\bb{\mathbf{f}}+1},
(\mathfrak{s}^{(2)}_{i})_{i\in\bb{\mathbf{f}}},
(W^{(1)}_{i})_{i\in\bb{\mathbf{f}}}
\right).
$$

Let us consider $\mathfrak{s}^{(2)}_{0}$, the first second-order rewrite rule appearing in $\mathfrak{Q}'^{(2)}$, and let us assume that it has the form 
$([K'_{0}]_{f_{0}}
,[L'_{0}]_{f_{0}})$ 
for some suitable pair of path terms classes $[K'_{0}]_{f_{0}}$ and $[L'_{0}]_{f_{0}}$ in $[\mathrm{PT}_{\boldsymbol{\mathcal{A}}}]_{f_{0}}$. Let $W^{(1)}_{0}$ be the initial first-order translation occurring in $\mathfrak{Q}'^{(2)}$.

Since $\mathfrak{Q}'^{(2)}$ is a second-order path, we have that 
\begin{itemize}
\item[(i)]
$
W^{(1)}_{0}(
K'_{0}
)
\in [S_{0}]_{s};
$
\item[(ii)]
$
W^{(1)}_{0}(
L'_{0}
)
\in [S_{1}]_{s}.
$
\end{itemize}

Let us recall that $W^{(1)}_{0}$ is a first-order translation of type $\tau$. Thus there will exists an index $j_{4}\in \bb{\mathbf{s}}$,  a family of paths $(\mathfrak{Q}'_{k})_{k\in j_{4}}\in \prod_{k\in j_{4}}\mathrm{Pth}_{\boldsymbol{\mathcal{A}},s_{k}}$, a family of paths $(\mathfrak{Q}'_{l})_{l\in \bb{\mathbf{s}}-(j_{4}+1)}\in \prod_{l\in \bb{\mathbf{s}}-(j_{4}+1)}\mathrm{Pth}_{\boldsymbol{\mathcal{A}},s_{l}}$ and a first-order translation $W^{(1)'}_{0}\in\mathrm{Tl}_{f_{0}}(\mathbf{T}_{\Sigma^{\boldsymbol{\mathcal{A}}}}(X))_{s_{j_{4}}}$ such that 
\begin{multline*}
W^{(1)}_{0}=\tau^{\mathbf{T}_{\Sigma^{\boldsymbol{\mathcal{A}}}}(X)}\left(
\mathrm{CH}^{(1)}_{s_{0}}\left(\mathfrak{Q}'_{0}\right),
\cdots,
\mathrm{CH}^{(1)}_{s_{j_{4}-1}}\left(\mathfrak{Q}'_{j_{4}-1}\right),
\right.
\\
\left.
W^{(1)'}_{0},
\mathrm{CH}^{(1)}_{s_{j_{4}+1}}\left(\mathfrak{Q}'_{j_{4}+1}\right),
\cdots,
\mathrm{CH}^{(1)}_{s_{\bb{\mathbf{s}}-1}}\left(\mathfrak{Q}'_{\bb{\mathbf{s}}-1}\right)
\right).
\end{multline*}

Moreover, the following chain of equalities holds
\begin{flushleft}
$\left[
S_{0}
\right]_{s}$
\allowdisplaybreaks
\begin{align*}
\,
&=
\left[
W^{(1)}_{0}\left(
K'_{0}
\right)
\right]_{s}
\tag{1}
\\&=
\biggr[
\tau^{\mathbf{PT}_{\boldsymbol{\mathcal{A}}}}
\biggr(
\mathrm{CH}^{(1)}_{s_{0}}\left(\mathfrak{Q}'_{0}\right),
\cdots,
\mathrm{CH}^{(1)}_{s_{j_{4}-1}}\left(\mathfrak{Q}'_{j_{4}-1}\right),
\\&\qquad\qquad\qquad
W^{(1)'}_{0}\left(
K'_{0}
\right),
\mathrm{CH}^{(1)}_{s_{j_{4}+1}}\left(\mathfrak{Q}'_{j_{4}+1}\right),
\cdots,
\mathrm{CH}^{(1)}_{s_{\bb{\mathbf{s}}-1}}\left(\mathfrak{Q}'_{\bb{\mathbf{s}}-1}\right)
\biggr)
\biggr]_{s}
\tag{2}
\\&=
\tau^{
[\mathbf{PT}_{\boldsymbol{\mathcal{A}}}]
}
\Biggr(
\left[
\mathrm{CH}^{(1)}_{s_{0}}\left(\mathfrak{Q}'_{0}\right)
\right]_{s_{0}},
\cdots,
\left[
\mathrm{CH}^{(1)}_{s_{j_{4}-1}}\left(\mathfrak{!}'_{j_{4}-1}\right)
\right]_{s_{j_{4}-1}},
\\&\qquad\qquad
\left[
W^{(1)'}_{0}\left(
K'_{0}
\right)
\right]_{s_{j_{4}}},
\\&\qquad\qquad\qquad\quad
\left[
\mathrm{CH}^{(1)}_{s_{j_{4}+1}}\left(\mathfrak{Q}'_{j_{4}+1}\right)
\right]_{s_{j_{4}+1}},
\cdots,
\left[
\mathrm{CH}^{(1)}_{s_{\bb{\mathbf{s}}-1}}\left(\mathfrak{Q}'_{\bb{\mathbf{s}}-1}\right)
\right]_{s_{\bb{\mathbf{s}}-1}}
\Biggr).
\tag{3}
\end{align*}
\end{flushleft}

The first equality follows from item (ii) above; the second equality unravels the description of the initial first-order translation occurring in $\mathfrak{Q}'^{(2)}$; finally, the last equality recovers the interpretation of the operation symbol $\tau$ in the many-sorted partial $\Sigma^{\boldsymbol{\mathcal{A}}}$-algebra $[\mathbf{PT}_{\boldsymbol{\mathcal{A}}}]$, as described in Proposition~\ref{PPTQCatAlg}. 

By Lemma~\ref{LDPthExtract}, we have the following family of equalities
\allowdisplaybreaks
\begin{align*}
\mathrm{sc}^{([1],2)}_{s_{0}}\left(
\mathfrak{Q}'^{(2)}_{0}
\right)
&=
\left[
\mathrm{CH}^{(1)}_{s_{0}}\left(\mathfrak{Q}'_{0}\right)
\right]_{s_{0}};
\\
\vdots\qquad&=\qquad\vdots
\\
\mathrm{sc}^{([1],2)}_{s_{j_{3}}}
\left(
\mathfrak{Q}'^{(2)}_{j_{3}}
\right)
&=
\left[
W^{(1)'}_{0}\left(
K'_{0}\right)
\right]_{s_{j_{4}}};
\\
\vdots\qquad&=\qquad\vdots
\\
\mathrm{sc}^{([1],2)}_{s_{\bb{\mathbf{s}}-1}}
\left(
\mathfrak{Q}'^{(2)}_{s_{\bb{\mathbf{s}}-1}}
\right)
&=
\left[
\mathrm{CH}^{(1)}_{s_{\bb{\mathbf{s}}-1}}\left(\mathfrak{Q}'_{\bb{\mathbf{s}}-1}\right)
\right]_{s_{\bb{\mathbf{s}}-1}}.
\end{align*}

Since we have already proved that, for every $j\in\bb{\mathbf{s}}$, the pair $(\mathfrak{P}'^{(2)}_{j}, \mathfrak{Q}'^{(2)}_{j})$ is in $\mathrm{Ker}(\mathrm{CH}^{(2)})_{s_{j}}$, then we have, by Lemma~\ref{LDCH}, that, for every $j\in\bb{\mathbf{s}}$, it happens that 
$$
\mathrm{sc}^{([1],2)}_{s_{j}}\left(
\mathfrak{P}^{(2)}_{j}
\right)
=
\mathrm{sc}^{([1],2)}_{s_{j}}\left(
\mathfrak{Q}^{(2)}_{j}
\right).
$$

Therefore, we have the following family of equalities between path term classes
\allowdisplaybreaks
\begin{align*}
\left[
\mathrm{CH}^{(1)}_{s_{0}}\left(\mathfrak{P}'_{0}\right)
\right]_{s_{0}}
&=
\left[
\mathrm{CH}^{(1)}_{s_{0}}\left(\mathfrak{Q}'_{0}\right)
\right]_{s_{0}};
\\
&\vdots
\\
\left[
V^{(1)'}_{0}\left(
M'_{0}
\right)
\right]_{s_{j_{3}}}
&=
\left[
\mathrm{CH}^{(1)}_{s_{j_{3}}}\left(\mathfrak{Q}'_{j_{3}}\right)
\right]_{s_{j_{3}}};
\\
&\vdots
\\
\left[
\mathrm{CH}^{(1)}_{s_{j_{4}}}\left(\mathfrak{P}'_{j_{4}}\right)
\right]_{s_{j_{4}}}
&=
\left[
W^{(1)'}_{0}\left(
K'_{0}
\right)
\right]_{s_{j_{4}}};
\\
&\vdots
\\
\left[
\mathrm{CH}^{(1)}_{s_{\bb{\mathbf{s}}-1}}\left(\mathfrak{P}'_{\bb{\mathbf{s}}-1}\right)
\right]_{s_{\bb{\mathbf{s}}-1}}
&=
\left[
\mathrm{CH}^{(1)}_{s_{\bb{\mathbf{s}}-1}}\left(\mathfrak{Q}'_{\bb{\mathbf{s}}-1}\right)
\right]_{s_{\bb{\mathbf{s}}-1}}.
\end{align*}

Now, it is time to put all the above information together. According to item (ii) in the hypothesis, $\mathfrak{P}'^{(2)}\circ_{s}^{1\mathbf{Pth}_{\boldsymbol{\mathcal{A}}^{(2)}}}\mathfrak{P}^{(2)}$ is a head-constant echelonless second-order path of length at least one that is not coherent and $i=\bb{\mathfrak{P}^{(2)}}-1=\bb{\mathbf{c}}-1$ is the greatest index for which $(\mathfrak{P}'^{(2)}\circ_{s}^{1\mathbf{Pth}_{\boldsymbol{\mathcal{A}}^{(2)}}}\mathfrak{P}^{(2)})^{0,i}$ is coherent. Hence, $(\mathfrak{P}'^{(2)}\circ_{s}^{1\mathbf{Pth}_{\boldsymbol{\mathcal{A}}^{(2)}}}\mathfrak{P}^{(2)})^{0,i+1}$ is no longer coherent. Let us recall that $i+1=\bb{\mathbf{c}}$, and the $i+1$-th translation of $\mathfrak{P}'^{(2)}\circ_{s}^{1\mathbf{Pth}_{\boldsymbol{\mathcal{A}}^{(2)}}}\mathfrak{P}^{(2)}$ is given by $V^{(1)}_{0}$, the initial first-order translation of $\mathfrak{P}'^{(2)}$.

Taking into account that $\mathrm{sc}^{([1],2)}_{s}(\mathfrak{P}'^{(2)})=\mathrm{tg}^{([1],2)}_{s}(\mathfrak{P}^{(2)})$, we have that the following equality holds
\allowdisplaybreaks
\begin{multline*}
\tau^{
[\mathbf{PT}_{\boldsymbol{\mathcal{A}}}]
}
\Biggr(
\left[
\mathrm{CH}^{(1)}_{s_{0}}\left(\mathfrak{P}'_{0}\right)
\right]_{s_{0}},
\cdots,
\left[
\mathrm{CH}^{(1)}_{s_{j_{3}-1}}\left(\mathfrak{P}'_{j_{3}-1}\right)
\right]_{s_{j_{3}-1}},
\left[
V^{(1)'}_{0}\left(
M'_{0}
\right)
\right]_{s_{j_{3}}},
\\
\left[
\mathrm{CH}^{(1)}_{s_{j_{3}+1}}\left(\mathfrak{P}'_{j_{3}+1}\right)
\right]_{s_{j_{3}+1}},
\cdots,
\left[
\mathrm{CH}^{(1)}_{s_{\bb{\mathbf{s}}-1}}\left(\mathfrak{P}'_{\bb{\mathbf{s}}-1}\right)
\right]_{s_{\bb{\mathbf{s}}-1}}
\Biggr)
\\=
\tau^{
[\mathbf{PT}_{\boldsymbol{\mathcal{A}}}]
}
\Biggr(
\left[
\mathrm{CH}^{(1)}_{s_{0}}\left(\mathfrak{P}_{0}\right)
\right]_{s_{0}},
\cdots,
\left[
\mathrm{CH}^{(1)}_{s_{j_{1}-1}}\left(\mathfrak{P}_{j_{1}-1}\right)
\right]_{s_{j_{1}-1}},
\left[
T^{(1)'}_{\bb{\mathbf{c}}-1}\left(
N_{\bb{\mathbf{c}}-1}
\right)
\right]_{s_{j_{1}}},
\\
\left[
\mathrm{CH}^{(1)}_{s_{j_{1}+1}}\left(\mathfrak{P}_{j_{1}+1}\right)
\right]_{s_{j_{1}+1}},
\cdots,
\left[
\mathrm{CH}^{(1)}_{s_{\bb{\mathbf{s}}-1}}\left(\mathfrak{P}_{\bb{\mathbf{s}}-1}\right)
\right]_{s_{\bb{\mathbf{s}}-1}}
\Biggr).
\end{multline*}

Since the coherence breaks at position $i+1$, from the above equation we  cannot infer that all the appearing path term classes constituting the overall class are equal, i.e., there exists at least an index $j_{5}\in\bb{\mathbf{s}}$ for which the $j_{5}$-th path term class on the left differs from the $j_{5}$-th path term class on the right.

Now, for the case of the interaction between the second-order paths $\mathfrak{Q}^{(2)}$ and $\mathfrak{Q}'^{(2)}$. Since $\mathrm{sc}^{([1],2)}_{s}(\mathfrak{Q}'^{(2)})=\mathrm{tg}^{([1],2)}_{s}(\mathfrak{Q}^{(2)})$, the following equality holds
\allowdisplaybreaks
\begin{multline*}
\tau^{
[\mathbf{PT}_{\boldsymbol{\mathcal{A}}}]
}
\Biggr(
\left[
\mathrm{CH}^{(1)}_{s_{0}}\left(\mathfrak{Q}'_{0}\right)
\right]_{s_{0}},
\cdots,
\left[
\mathrm{CH}^{(1)}_{s_{j_{4}-1}}\left(\mathfrak{!}'_{j_{4}-1}\right)
\right]_{s_{j_{4}-1}},
\left[
W^{(1)'}_{0}\left(
K'_{0}
\right)
\right]_{s_{j_{4}}},
\\
\left[
\mathrm{CH}^{(1)}_{s_{j_{4}+1}}\left(\mathfrak{Q}'_{j_{4}+1}\right)
\right]_{s_{j_{4}+1}},
\cdots,
\left[
\mathrm{CH}^{(1)}_{s_{\bb{\mathbf{s}}-1}}\left(\mathfrak{Q}'_{\bb{\mathbf{s}}-1}\right)
\right]_{s_{\bb{\mathbf{s}}-1}}
\Biggr)
\\=
\tau^{
[\mathbf{PT}_{\boldsymbol{\mathcal{A}}}]
}
\Biggr(
\left[
\mathrm{CH}^{(1)}_{s_{0}}\left(\mathfrak{Q}_{0}\right)
\right]_{s_{0}},
\cdots,
\left[
\mathrm{CH}^{(1)}_{s_{j_{2}-1}}\left(\mathfrak{Q}_{j_{2}-1}\right)
\right]_{s_{j_{2}-1}},
\left[
U^{(1)'}_{\bb{\mathbf{d}}-1}\left(
L_{\bb{\mathbf{d}}-1}
\right)
\right]_{s_{j_{2}}},
\\
\left[
\mathrm{CH}^{(1)}_{s_{j_{2}+1}}\left(\mathfrak{Q}_{j_{2}+1}\right)
\right]_{s_{j_{2}+1}},
\cdots,
\left[
\mathrm{CH}^{(1)}_{s_{\bb{\mathbf{s}}-1}}\left(\mathfrak{Q}_{\bb{\mathbf{s}}-1}\right)
\right]_{s_{\bb{\mathbf{s}}-1}}
\Biggr).
\end{multline*}

Taking into account the previous equalities, from the above equation we cannot infer that all the appearing path term classes constituting the overall class are equal, i.e., for the index $j_{5}\in\bb{\mathbf{s}}$ previously encountered, the $j_{5}$-th path term class on the left differs from the $j_{5}$-th path term class on the right.

Therefore, $\mathfrak{Q}'^{(2)}\circ_{s}^{1\mathbf{Pth}_{\boldsymbol{\mathcal{A}}^{(2)}}}\mathfrak{Q}^{(2)}$ is a head-constant echelonless second-order path that is not coherent. Moreover, $i=\bb{\mathfrak{P}^{(2)}}-1=\bb{\mathfrak{Q}^{(2)}}-1$ is the greatest index in $\bb{\mathfrak{Q}'^{(2)}\circ_{s}^{1\mathbf{Pth}_{\boldsymbol{\mathcal{A}}^{(2)}}}\mathfrak{Q}^{(2)}}$ for which $(\mathfrak{Q}'^{(2)}\circ_{s}^{1\mathbf{Pth}_{\boldsymbol{\mathcal{A}}^{(2)}}}\mathfrak{Q}^{(2)})^{0,i}$ is coherent.

This completes the proof.
\end{proof}

\section{
\texorpdfstring
{A  congruence on second-order paths}
{A  congruence on second-order paths}
}
We are now in position to prove that the kernel of the second-order Curry-Howard mapping is a $\Sigma^{\boldsymbol{\mathcal{A}}^{(2)}}$-congruence on $\mathbf{Pth}_{\boldsymbol{\mathcal{A}}^{(2)}}$.

\begin{restatable}{proposition}{PDCHCong}
\label{PDCHCong}
$\mathrm{Ker}(\mathrm{CH}^{(2)})$ is a closed $\Sigma^{\boldsymbol{\mathcal{A}}^{(2)}}$-congruence on $\mathbf{Pth}_{\boldsymbol{\mathcal{A}}^{(2)}}$.
\end{restatable}

\begin{proof}
By Proposition~\ref{PDCHHom}, the second-order Curry-Howard mapping is a $\Sigma$-homomorphism from $\mathbf{Pth}_{\boldsymbol{\mathcal{A}}^{(2)}}^{(0,2)}$, the $\Sigma$-reduct of the $\Sigma^{\boldsymbol{\mathcal{A}}^{(2)}}$-algebra $\mathbf{Pth}_{\boldsymbol{\mathcal{A}}^{(2)}}$, to $\mathbf{T}_{\Sigma^{\boldsymbol{\mathcal{A}}^{(2)}}}^{(0,2)}(X)$, the $\Sigma$-reduct of the free $\Sigma^{\boldsymbol{\mathcal{A}}^{(2)}}$-algebra $\mathbf{T}_{\Sigma^{\boldsymbol{\mathcal{A}}^{(2)}}}(X)$. Therefore $\mathrm{Ker}(\mathrm{CH}^{(2)})$ is a $\Sigma$-congruence.

Thus, we only need to check the compatibility of $\mathrm{Ker}(\mathrm{CH}^{(2)})$ with the categorial operations. The addition of new constants to the signature, obviously, does not matter. Let us begin with the $0$-source operation symbol.

\begin{claim}\label{CDCHScZ}
Let $s$ be a sort in $S$ and let $\mathrm{sc}^{0}_{s}$ be the $0$-source operation symbol of type $s$ in $\Sigma^{\boldsymbol{\mathcal{A}}^{(2)}}_{s,s}$. Let $\mathfrak{P}^{(2)}$ and $\mathfrak{Q}^{(2)}$ be two second-order paths in $\mathrm{Pth}_{\boldsymbol{\mathcal{A}}^{(2)},s}$ such that $(\mathfrak{P}^{(2)},\mathfrak{Q}^{(2)})\in\mathrm{Ker}(\mathrm{CH}^{(2)})_{s}$, then 
$$
\left(
\mathrm{sc}^{0\mathbf{Pth}_{\boldsymbol{\mathcal{A}}^{(2)}}}_{s}\left(
\mathfrak{P}^{(2)}
\right), \mathrm{sc}^{0\mathbf{Pth}_{\boldsymbol{\mathcal{A}}^{(2)}}}_{s}
\left(
\mathfrak{Q}^{(2)}
\right)
\right)\in\mathrm{Ker}\left(\mathrm{CH}^{(2)}
\right)_{s}.$$
\end{claim}
 
To prove the claim, note that according to Proposition~\ref{CDCH}, we have that 
\[
\mathrm{sc}^{(0,2)}_{s}\left(\mathfrak{P}^{(2)}\right)=
\mathrm{sc}^{(0,2)}_{s}\left(\mathfrak{Q}^{(2)}\right).
\]

Following the interpretation of the $0$-source operation symbol in the many-sorted partial $\Sigma^{\boldsymbol{\mathcal{A}}^{(2)}}$-algebra $\mathbf{Pth}_{\boldsymbol{\mathcal{A}}^{(2)}}$ introduced in Proposition~\ref{PDPthCatAlg}, we conclude that 
\allowdisplaybreaks
\begin{align*}
\mathrm{sc}^{0\mathbf{Pth}_{\boldsymbol{\mathcal{A}}^{(2)}}}_{s}\left(\mathfrak{P}^{(2)}\right)
&=
\mathrm{ip}^{(2,0)\sharp}_{s}\left(
\mathrm{sc}^{(0,2)}_{s}\left(
\mathfrak{P}^{(2)}
\right)\right)
\\&=
\mathrm{ip}^{(2,0)\sharp}_{s}\left(
\mathrm{sc}^{(0,2)}_{s}\left(
\mathfrak{Q}^{(2)}
\right)\right)
\\&=
\mathrm{sc}^{0\mathbf{Pth}_{\boldsymbol{\mathcal{A}}^{(2)}}}_{s}\left(
\mathfrak{Q}^{(2)}\right).
\end{align*}

Claim~\ref{CDCHScZ} follows by reflexivity. This finishes the proof of Claim~\ref{CDCHScZ}. Therefore, $\mathrm{Ker}(\mathrm{CH}^{(2)})$ is compatible with the $0$-source operation symbol.

 Let us now consider the $0$-target operation symbol.
\begin{claim}\label{CDCHTgZ}
Let $s$ be a sort in $S$ and let $\mathrm{tg}^{0}_{s}$ be the $0$-target operation symbol of type $s$ in $\Sigma^{\boldsymbol{\mathcal{A}}^{(2)}}_{s,s}$. Let $\mathfrak{P}^{(2)}$ and $\mathfrak{Q}^{(2)}$ be two second-order paths in $\mathrm{Pth}_{\boldsymbol{\mathcal{A}}^{(2)},s}$ such that $(\mathfrak{P}^{(2)},\mathfrak{Q}^{(2)})\in\mathrm{Ker}(\mathrm{CH}^{(2)})_{s}$, then 
$$
\left(
\mathrm{tg}^{0\mathbf{Pth}_{\boldsymbol{\mathcal{A}}^{(2)}}}_{s}\left(
\mathfrak{P}^{(2)}
\right), \mathrm{tg}^{0\mathbf{Pth}_{\boldsymbol{\mathcal{A}}^{(2)}}}_{s}\left(
\mathfrak{Q}^{(2)}
\right)\right)\in\mathrm{Ker}\left(\mathrm{CH}^{(2)}\right)_{s}.$$
\end{claim}

The proof of Claim~\ref{CDCHTgZ} follows from a similar argument to that presented in the proof of Claim~\ref{CDCHScZ}.

Let us now consider the $0$-composition operation symbol.
\begin{claim}\label{CDCHCompZ}
Let $s$ be a sort in $S$ and let $\circ^{0}_{s}$ be the $0$-composition operation symbol of type $s$ in $\Sigma^{\boldsymbol{\mathcal{A}}^{(2)}}_{ss,s}$. Let $\mathfrak{P}'^{(2)},\mathfrak{P}^{(2)}, \mathfrak{Q}'^{(2)}$ and $\mathfrak{Q}^{(2)}$ be four second-order paths in $\mathrm{Pth}_{\boldsymbol{\mathcal{A}}^{(2)},s}$ such that $(\mathfrak{P}'^{(2)},\mathfrak{Q}'^{(2)})\in\mathrm{Ker}(\mathrm{CH}^{(2)})_{s}$ and $(\mathfrak{P}^{(2)},\mathfrak{Q}^{(2)})\in\mathrm{Ker}(\mathrm{CH}^{(2)})_{s}$ then, when defined, 
$$\left(
\mathfrak{P}'^{(2)}
\circ^{0\mathbf{Pth}_{\boldsymbol{\mathcal{A}}^{(2)}}}_{s}
\mathfrak{P}^{(2)}, 
\mathfrak{Q}'^{(2)}
\circ^{0\mathbf{Pth}_{\boldsymbol{\mathcal{A}}^{(2)}}}_{s}
\mathfrak{Q}^{(2)}
\right)
\in\mathrm{Ker}(\mathrm{CH}^{(2)})_{s}.$$
\end{claim}

Let us assume that the operation $\mathfrak{P}'^{(2)}
\circ^{0\mathbf{Pth}_{\boldsymbol{\mathcal{A}}^{(2)}}}_{s}
\mathfrak{P}^{(2)}$ is defined. But, from Proposition~\ref{PDPthCatAlg}, this is the case when
$
\mathrm{sc}^{(0,2)}_{s}(
\mathfrak{P}'^{(2)}
)
=
\mathrm{tg}^{(0,2)}_{s}(
\mathfrak{P}^{(2)}
).
$
Since $(\mathfrak{P}'^{(2)},\mathfrak{Q}'^{(2)})\in\mathrm{Ker}(\mathrm{CH}^{(2)})_{s}$ and $(\mathfrak{P}^{(2)},\mathfrak{Q}^{(2)})\in\mathrm{Ker}(\mathrm{CH}^{(2)})_{s}$ we have, by Corollary~\ref{CDCH}, that 
\begin{align*}
\mathrm{sc}^{(0,2)}_{s}\left(
\mathfrak{P}'^{(2)}
\right)
&=
\mathrm{sc}^{(0,2)}_{s}\left(
\mathfrak{Q}'^{(2)}
\right);
&
\mathrm{tg}^{(0,2)}_{s}\left(
\mathfrak{P}^{(2)}
\right)
&=
\mathrm{tg}^{(0,2)}_{s}\left(
\mathfrak{Q}^{(2)}
\right).
\end{align*}

All in all, we can affirm that 
$
\mathrm{sc}^{(0,2)}_{s}(
\mathfrak{Q}'^{(2)}
)
=
\mathrm{tg}^{(0,2)}_{s}(
\mathfrak{Q}^{(2)}
).
$
Therefore, according to Proposition~\ref{PDPthCatAlg}, the $0$-composition $\mathfrak{Q}'^{(2)}
\circ^{0\mathbf{Pth}_{\boldsymbol{\mathcal{A}}^{(2)}}}_{s}
\mathfrak{Q}^{(2)}$ is defined

The other implication follows by the symmetry in the argument.
Consequently, $\mathfrak{P}'^{(2)}
\circ^{0\mathbf{Pth}_{\boldsymbol{\mathcal{A}}^{(2)}}}_{s}
\mathfrak{P}^{(2)}$ is defined if, and only if,
$\mathfrak{Q}'^{(2)}
\circ^{0\mathbf{Pth}_{\boldsymbol{\mathcal{A}}^{(2)}}}_{s}
\mathfrak{Q}^{(2)}$ is defined.

To prove the claim, we will distinguish two cases according to the nature of the pair $\mathfrak{P}'^{(2)}$ and $\mathfrak{P}^{(2)}$. It could be the case that either (1) both $\mathfrak{P}'^{(2)}$ and $\mathfrak{P}^{(2)}$ are $(2,[1])$-identity second-order paths or (2) at least one of the second-order paths $\mathfrak{P}'^{(2)}$ and $\mathfrak{P}^{(2)}$ has length at least one.

If (1), i.e., if $\mathfrak{P}'^{(2)}$ and $\mathfrak{P}^{(2)}$ are $(2,[1])$-identity second-order paths then, since $(\mathfrak{P}'^{(2)},\mathfrak{Q}'^{(2)})$ and $(\mathfrak{P}^{(2)},\mathfrak{Q}^{(2)})$ are pairs in $\mathrm{Ker}(\mathrm{CH}^{(2)})_{s}$, then, by Corollary~\ref{CDCHUId}, we have that $\mathfrak{P}'^{(2)}=\mathfrak{Q}'^{(2)}$ and $\mathfrak{P}^{(2)}=\mathfrak{Q}^{(2)}$. Consequently,
$$
\mathrm{CH}^{(2)}_{s}\left(
\mathfrak{P}'^{(2)}
\circ^{0\mathbf{Pth}_{\boldsymbol{\mathcal{A}}^{(2)}}}_{s}
\mathfrak{P}^{(2)}
\right)
=
\mathrm{CH}^{(2)}_{s}\left(
\mathfrak{Q}'^{(2)}
\circ^{0\mathbf{Pth}_{\boldsymbol{\mathcal{A}}^{(2)}}}_{s}
\mathfrak{Q}^{(2)}
\right).
$$

The case in which both $\mathfrak{P}'^{(2)}$ and $\mathfrak{P}^{(2)}$ are $(2,[1])$-identity second-order paths follows.

If~(2), i.e., if at least one of the second-order paths $\mathfrak{P}'^{(2)}$ and $\mathfrak{P}^{(2)}$ has length at least one, we can assume, without loss of generality, that it is the case that $\mathfrak{P}^{(2)}$ has length at least one. Then, since $(\mathfrak{P}^{(2)},\mathfrak{Q}^{(2)})\in\mathrm{Ker}(\mathrm{CH}^{(2)})_{s}$, we have according to Corollary~\ref{CDCHUId} that $\mathfrak{Q}^{(2)}$ is also a second-order path of length at least one.

Then, by Corollary~\ref{CDPthWB}, we have that the second-order paths 
$\mathfrak{P}'^{(2)}
\circ^{0\mathbf{Pth}_{\boldsymbol{\mathcal{A}}^{(2)}}}_{s}
\mathfrak{P}^{(2)}$ and 
$\mathfrak{Q}'^{(2)}
\circ^{0\mathbf{Pth}_{\boldsymbol{\mathcal{A}}^{(2)}}}_{s}
\mathfrak{Q}^{(2)}$ are both head-constant coherent and echelonless second-order paths. Moreover, by Proposition~\ref{PDRecov}, the second-order path extraction procedure from Lemma~\ref{LDPthExtract} applied to the second-order paths $\mathfrak{P}'^{(2)}
\circ^{0\mathbf{Pth}_{\boldsymbol{\mathcal{A}}^{(2)}}}_{s}
\mathfrak{P}^{(2)}$ and 
$\mathfrak{Q}'^{(2)}
\circ^{0\mathbf{Pth}_{\boldsymbol{\mathcal{A}}^{(2)}}}_{s}
\mathfrak{Q}^{(2)}$ retrieves, respectively, the pair of second-order paths $(\mathfrak{P}'^{(2)},\mathfrak{P}^{(2)})$ and the pair of second-order paths 
$(\mathfrak{Q}'^{(2)},\mathfrak{Q}^{(2)})$.

Thus, according to Definition~\ref{DDCH}, the value of the second-order Curry-Howard mapping at $\mathfrak{P}'^{(2)}
\circ^{0\mathbf{Pth}_{\boldsymbol{\mathcal{A}}^{(2)}}}_{s}
\mathfrak{P}^{(2)}$ and 
$\mathfrak{Q}'^{(2)}
\circ^{0\mathbf{Pth}_{\boldsymbol{\mathcal{A}}^{(2)}}}_{s}
\mathfrak{Q}^{(2)}$ is given, respectively, by the following terms
\begin{align*}
\mathrm{CH}^{(2)}_{s}\left(
\mathfrak{P}'^{(2)}
\circ^{0\mathbf{Pth}_{\boldsymbol{\mathcal{A}}^{(2)}}}_{s}
\mathfrak{P}^{(2)}
\right)&=
\mathrm{CH}^{(2)}_{s}\left(
\mathfrak{P}'^{(2)}
\right)
\circ^{0\mathbf{T}_{\Sigma^{\boldsymbol{\mathcal{A}}^{(2)}}}(X)}_{s}
\mathrm{CH}^{(2)}_{s}\left(
\mathfrak{P}^{(2)}
\right);
\\
\mathrm{CH}^{(2)}_{s}\left(
\mathfrak{Q}'^{(2)}
\circ^{0\mathbf{Pth}_{\boldsymbol{\mathcal{A}}^{(2)}}}_{s}
\mathfrak{Q}^{(2)}
\right)&=
\mathrm{CH}^{(2)}_{s}\left(
\mathfrak{Q}'^{(2)}
\right)
\circ^{0\mathbf{T}_{\Sigma^{\boldsymbol{\mathcal{A}}^{(2)}}}(X)}_{s}
\mathrm{CH}^{(2)}_{s}\left(
\mathfrak{Q}^{(2)}
\right).
\end{align*}

Since $(\mathfrak{P}'^{(2)},\mathfrak{Q}'^{(2)})$ and $(\mathfrak{P}^{(2)},\mathfrak{Q}^{(2)})$ are pairs in $\mathrm{Ker}(\mathrm{CH}^{(2)})_{s}$, it follows that 
$$
\left(
\mathfrak{P}'^{(2)}
\circ^{0\mathbf{Pth}_{\boldsymbol{\mathcal{A}}^{(2)}}}_{s}
\mathfrak{P}^{(2)},
\mathfrak{Q}'^{(2)}
\circ^{0\mathbf{Pth}_{\boldsymbol{\mathcal{A}}^{(2)}}}_{s}
\mathfrak{Q}^{(2)}
\right)
\in\mathrm{Ker}(\mathrm{CH}^{(2)})_{s}.
$$

This finishes the proof of Claim~\ref{CDCHCompZ}. Therefore, $\mathrm{Ker}(\mathrm{CH}^{(2)})_{s}$ is compatible with the $0$-composition operation symbol.

Regarding the $1$-source and the $1$-target operations, we note that, for every sort $s\in S$, if $(\mathfrak{P}^{(2)},\mathfrak{Q}^{(2)})$ is a pair of second-order paths in $\mathrm{Ker}(\mathrm{CH}^{(2)})_{s}$ then, by Lemma~\ref{LDCH}, we have that
\begin{align*}
\mathrm{sc}^{([1],2)}_{s}\left(
\mathfrak{P}^{(2)}
\right)&=
\mathrm{sc}^{([1],2)}_{s}\left(
\mathfrak{Q}^{(2)}
\right),
&
\mathrm{tg}^{([1],2)}_{s}\left(
\mathfrak{P}^{(2)}
\right)&=
\mathrm{tg}^{([1],2)}_{s}\left(
\mathfrak{Q}^{(2)}
\right).
\end{align*}

Let us recall from Proposition~\ref{PDPthDCatAlg} the interpretation of the operation symbols $\mathrm{sc}^{1}_{s}$ and $\mathrm{tg}^{1}_{s}$ in the partial many-sorted $\Sigma^{\boldsymbol{\mathcal{A}}^{(2)}}$-algebra $\mathbf{Pth}_{\boldsymbol{\mathcal{A}}^{(2)}}$. Hence we have that
\allowdisplaybreaks
\begin{align*}
\mathrm{sc}^{1\mathbf{Pth}_{\boldsymbol{\mathcal{A}}^{(2)}}}_{s}\left(
\mathfrak{P}^{(2)}
\right)
&=
\mathrm{ip}^{(2,[1])\sharp}_{s}\left(
\mathrm{sc}^{([1],2)}_{s}\left(
\mathfrak{P}^{(2)}
\right)\right)
\\&=
\mathrm{ip}^{(2,[1])\sharp}_{s}\left(
\mathrm{sc}^{([1],2)}_{s}\left(
\mathfrak{Q}^{(2)}
\right)\right)
\\&=
\mathrm{sc}^{1\mathbf{Pth}_{\boldsymbol{\mathcal{A}}^{(2)}}}_{s}\left(
\mathfrak{Q}^{(2)}
\right);
\\
\\
\mathrm{tg}^{1\mathbf{Pth}_{\boldsymbol{\mathcal{A}}^{(2)}}}_{s}\left(
\mathfrak{P}^{(2)}
\right)
&=
\mathrm{ip}^{(2,[1])\sharp}_{s}\left(
\mathrm{tg}^{([1],2)}_{s}\left(
\mathfrak{P}^{(2)}
\right)\right)
\\&=
\mathrm{ip}^{(2,[1])\sharp}_{s}\left(
\mathrm{tg}^{([1],2)}_{s}\left(
\mathfrak{Q}^{(2)}
\right)\right)
\\&=
\mathrm{tg}^{1\mathbf{Pth}_{\boldsymbol{\mathcal{A}}^{(2)}}}_{s}\left(
\mathfrak{Q}^{(2)}
\right).
\end{align*}

Since $\mathrm{sc}^{1\mathbf{Pth}_{\boldsymbol{\mathcal{A}}^{(2)}}}_{s}(
\mathfrak{P}^{(2)}
)=\mathrm{sc}^{1\mathbf{Pth}_{\boldsymbol{\mathcal{A}}^{(2)}}}_{s}(
\mathfrak{Q}^{(2)}
)$ and  $\mathrm{tg}^{1\mathbf{Pth}_{\boldsymbol{\mathcal{A}}^{(2)}}}_{s}(
\mathfrak{P}^{(2)}
)=\mathrm{tg}^{1\mathbf{Pth}_{\boldsymbol{\mathcal{A}}^{(2)}}}_{s}(
\mathfrak{Q}^{(2)}
)$, it follows that they have, respectively, the same image under the second-order Curry-Howard mapping.

Therefore, we are only left with the compatibility of $\mathrm{Ker}(\mathrm{CH}^{(2)})$ with the $1$-composition operator.

Let $s$ be a sort in $S$ and let $(\mathfrak{P}^{(2)},\mathfrak{Q}^{(2)})$ and $(\mathfrak{P}'^{(2)},\mathfrak{Q}'^{(2)})$ be two pairs of second-order paths in $\mathrm{Ker}(\mathrm{CH}^{(2)})_{s}$. Let us assume that the $1$-composition $\mathfrak{P}'^{(2)}\circ^{1\mathbf{Pth}_{\boldsymbol{\mathcal{A}}^{(2)}}}_{s}\mathfrak{P}^{(2)}$ is defined. Then, by Lemma~\ref{LDCH}, the $1$-composition $\mathfrak{Q}'^{(2)}\circ^{1\mathbf{Pth}_{\boldsymbol{\mathcal{A}}^{(2)}}}_{s}\mathfrak{Q}^{(2)}$ is also defined. We want to prove that their respective $1$-composition, i.e.,
$\mathfrak{P}'^{(2)}
\circ^{1\mathbf{Pth}_{\boldsymbol{\mathcal{A}}^{(2)}}}_{s}
\mathfrak{P}^{(2)}$ and $
\mathfrak{Q}'^{(2)}
\circ^{1\mathbf{Pth}_{\boldsymbol{\mathcal{A}}^{(2)}}}_{s}
\mathfrak{Q}^{(2)}
$ have the same value for the second-order Curry-Howard mapping.

We first study the case in which one of the second-order paths is a $(2,[1])$-identity second-order path.

\begin{claim}\label{CDUCompIp} Let $s$ be a sort in $S$ and let $(\mathfrak{P}^{(2)},\mathfrak{Q}^{(2)})$ and $(\mathfrak{P}'^{(2)},\mathfrak{Q}'^{(2)})$ be two pairs of second-order paths in $\mathrm{Ker}(\mathrm{CH}^{(2)})_{s}$. If $\mathfrak{P}^{(2)},\mathfrak{P}'^{(2)},\mathfrak{Q}^{(2)}$ or $\mathfrak{Q}'^{(2)}$ is a $(2,[1])$-identity second-order path then, when defined,
$$
\left(\mathfrak{P}'^{(2)}\circ^{1\mathbf{Pth}_{\boldsymbol{\mathcal{A}}^{(2)}}}_{s}\mathfrak{P}^{(2)}
,
\mathfrak{Q}'^{(2)}\circ^{1\mathbf{Pth}_{\boldsymbol{\mathcal{A}}^{(2)}}}_{s}\mathfrak{Q}^{(2)}
\right)\in\mathrm{Ker}\left(\mathrm{CH}^{(2)}\right)_{s}.
$$
\end{claim}

Assume that $\mathfrak{P}^{(2)}$ is a $(2,[1])$-identity second-order path.

Since 
$\mathfrak{P}'^{(2)}
\circ^{1\mathbf{Pth}_{\boldsymbol{\mathcal{A}}^{(2)}}}_{s}
\mathfrak{P}^{(2)}$
is defined, we have that $\mathrm{sc}^{([1],2)}_{s}(\mathfrak{P}'^{(2)})
=\mathrm{tg}^{([1],2)}_{s}(\mathfrak{P}^{(2)})$. Hence $\mathfrak{P}^{(2)}$ is the $(2,1)$-identity second-order path on the $([1],2)$-source of $\mathfrak{P}^{(2)}$, i.e., 
$\mathfrak{P}^{(2)}=
\mathrm{ip}^{(2,[1])\sharp}_{s}(
\mathrm{sc}^{([1],2)}_{s}(
\mathfrak{P}'^{(2)}))$. Hence the $1$-composition $\mathfrak{P}'^{(2)}\circ^{1\mathbf{Pth}_{\boldsymbol{\mathcal{A}}^{(2)}}}_{s}\mathfrak{P}^{(2)}$ reduces to $\mathfrak{P}'^{(2)}$.

On the other hand, since $(\mathfrak{P}^{(2)},\mathfrak{Q}^{(2)})$ is a pair in $\mathrm{Ker}(\mathrm{CH}^{(2)})_{s}$, we have, by Corollary~\ref{CDCHUId}, that $\mathfrak{P}^{(2)}=\mathfrak{Q}^{(2)}$. Hence, $\mathfrak{Q}^{(2)}=\mathrm{ip}^{(2,[1])\sharp}_{s}(\mathrm{sc}^{([1],2)}_{s}(\mathfrak{P}^{(2)}))$. Moreover, since $(\mathfrak{P}'^{(2)},\mathfrak{Q}'^{(2)})$ is a pair in $\mathrm{Ker}(\mathrm{CH}^{(2)})_{s}$, we have, by Lemma~\ref{LDCH}, that $\mathrm{sc}^{([1],2)}_{s}(\mathfrak{P}'^{(2)})=\mathrm{sc}^{([1],2)}_{s}(\mathfrak{Q}'^{(2)})$. Hence, the $1$-composition $\mathfrak{Q}'^{(2)}\circ^{1\mathbf{Pth}_{\boldsymbol{\mathcal{A}}^{(2)}}}_{s}\mathfrak{Q}^{(2)}$ reduces to $\mathfrak{Q}'^{(2)}$.

Then, we have that 
$$
\left(
\mathfrak{P}'^{(2)}
\circ^{1\mathbf{Pth}_{\boldsymbol{\mathcal{A}}^{(2)}}}_{s}
\mathfrak{P}^{(2)}
,
\mathfrak{Q}'^{(2)}
\circ^{1\mathbf{Pth}_{\boldsymbol{\mathcal{A}}^{(2)}}}_{s}
\mathfrak{Q}^{(2)}
\right)\in\mathrm{Ker}\left(\mathrm{CH}^{(2)}
\right)_{s}.
$$

The remaining cases are argued similarly. This finishes the proof of Claim~\ref{CDUCompIp}.

We are now in position to prove the general situation.

\begin{claim}\label{CDUCompProd} Let $s$ be a sort in $S$ and let $(\mathfrak{P}^{(2)},\mathfrak{Q}^{(2)})$ and $(\mathfrak{P}'^{(2)},\mathfrak{Q}'^{(2)})$ be two pairs of second-order paths in $\mathrm{Ker}(\mathrm{CH}^{(2)})_{s}$. When defined, 
$$
\left(\mathfrak{P}'^{(2)}\circ^{1\mathbf{Pth}_{\boldsymbol{\mathcal{A}}^{(2)}}}_{s}\mathfrak{P}^{(2)}
,
\mathfrak{Q}'^{(2)}\circ^{1\mathbf{Pth}_{\boldsymbol{\mathcal{A}}^{(2)}}}_{s}\mathfrak{Q}^{(2)}
\right)\in\mathrm{Ker}\left(\mathrm{CH}^{(2)}
\right)_{s}.
$$
\end{claim}

The proof is done by Artinian induction on $(\coprod\mathrm{Pth}_{\boldsymbol{\mathcal{A}}^{(2)}}, \leq_{\mathbf{Pth}_{\boldsymbol{\mathcal{A}}^{(2)}}})$.

\textsf{Base step of the Artinian induction.}

Let $(\mathfrak{P}'^{(2)}
\circ^{1\mathbf{Pth}_{\boldsymbol{\mathcal{A}}^{(2)}}}_{s}
\mathfrak{P}^{(2)},s)$ be a minimal element in $(\coprod\mathrm{Pth}_{\boldsymbol{\mathcal{A}}^{(2)}}, \leq_{\mathbf{Pth}_{\boldsymbol{\mathcal{A}}^{(2)}}})$. Then, by Proposition~\ref{PDMinimal}, the second-order path $\mathfrak{P}'^{(2)}
\circ^{1\mathbf{Pth}_{\boldsymbol{\mathcal{A}}^{(2)}}}_{s}
\mathfrak{P}^{(2)}$ is either a $(2,[1])$-identity second-order path or a   second-order echelon. In any case, either $\mathfrak{P}^{(2)}$ or $\mathfrak{P}'^{(2)}$ must be a $(2,[1])$-identity second-order path. The statement follows by Claim~\ref{CDUCompIp}.

\textsf{Inductive step of the Artinian induction.}

Let $(\mathfrak{P}'^{(2)}\circ^{1\mathbf{Pth}_{\boldsymbol{\mathcal{A}}^{(2)}}}_{s}\mathfrak{P}^{(2)},s)$ be a non-minimal element in $(\coprod\mathrm{Pth}_{\boldsymbol{\mathcal{A}}^{(2)}}, \leq_{\mathbf{Pth}_{\boldsymbol{\mathcal{A}}^{(2)}}})$. Let us suppose that, for every sort $t\in S$ and every second-order path $\mathfrak{P}'''^{(2)}\circ^{1\mathbf{Pth}_{\boldsymbol{\mathcal{A}}^{(2)}}}_{t}\mathfrak{P}''^{(2)}$ in $\mathrm{Pth}_{\boldsymbol{\mathcal{A}}^{(2)},t}$, if $(\mathfrak{P}'''^{(2)}\circ^{1\mathbf{Pth}_{\boldsymbol{\mathcal{A}}^{(2)}}}_{t}\mathfrak{P}''^{(2)},t)<_{\mathbf{Pth}_{\boldsymbol{\mathcal{A}}^{(2)}}}(\mathfrak{P}'^{(2)}\circ^{1\mathbf{Pth}_{\boldsymbol{\mathcal{A}}^{(2)}}}_{s}\mathfrak{P}^{(2)},s)$, then the statement holds for $\mathfrak{P}'''^{(2)}\circ^{1\mathbf{Pth}_{\boldsymbol{\mathcal{A}}^{(2)}}}_{t}\mathfrak{P}''^{(2)}$, i.e., for every pair of second-order paths $\mathfrak{Q}'''^{(2)}$ and $\mathfrak{Q}''^{(2)}$ in $\mathrm{Pth}_{\boldsymbol{\mathcal{A}}^{(2)},t}$, if $(\mathfrak{P}''^{(2)},\mathfrak{Q}''^{(2)})$ and $(\mathfrak{P}'''^{(2)},\mathfrak{Q}'''^{(2)})$ are pairs in $\mathrm{Ker}(\mathrm{CH}^{(2)})_{t}$ then
$$
\left(
\mathfrak{P}'''^{(2)}
\circ^{1\mathbf{Pth}_{\boldsymbol{\mathcal{A}}^{(2)}}}_{t}
\mathfrak{P}''^{(2)}
,
\mathfrak{Q}'''^{(2)}
\circ^{1\mathbf{Pth}_{\boldsymbol{\mathcal{A}}^{(2)}}}_{t}
\mathfrak{Q}''^{(2)}
\right)\in\mathrm{Ker}\left(\mathrm{CH}^{(2)}
\right)_{t}.
$$

Since $(\mathfrak{P}'^{(2)}\circ_{s}^{1\mathbf{Pth}_{\boldsymbol{\mathcal{A}}^{(2)}}}\mathfrak{P}^{(2)},s)$ is a non-minimal element in $(\coprod\mathrm{Pth}_{\boldsymbol{\mathcal{A}}^{(2)}}, \leq_{\mathbf{Pth}_{\boldsymbol{\mathcal{A}}^{(2)}}})$ and, by Claim~\ref{CDUCompIp}, we can assume that neither $\mathfrak{P}'^{(2)}$ nor $\mathfrak{P}^{(2)}$ are $(2,[1])$-identity second-order paths, we have, by Lemma~\ref{LDOrdI}, that $\mathfrak{P}'^{(2)}\circ^{1\mathbf{Pth}_{\boldsymbol{\mathcal{A}}^{(2)}}}\mathfrak{P}^{(2)}$ is either (1) a second-order path of length strictly greater than one containing at least one   second-order echelon or (2) an echelonless second-order path.

If~(1), then let $i\in\bb{\mathfrak{P}'^{(2)}\circ^{1\mathbf{Pth}_{\boldsymbol{\mathcal{A}}^{(2)}}}_{s}\mathfrak{P}^{(2)}}$ be the first index for which the one-step subpath $(\mathfrak{P}'^{(2)}\circ^{1\mathbf{Pth}_{\boldsymbol{\mathcal{A}}^{(2)}}}_{s}\mathfrak{P}^{(2)})^{i,i}$ of $\mathfrak{P}'^{(2)}\circ^{1\mathbf{Pth}_{\boldsymbol{\mathcal{A}}^{(2)}}}_{s}\mathfrak{P}^{(2)}$ is a   second-order echelon. We distinguish the cases~(1.1) $i=0$ and~(1.2) $i\neq 0$.

If~(1.1), i.e., $i=0$, since we are assuming that $\mathfrak{P}^{(2)}$ is not a $(2,[1])$-identity second-order path, we have that $\mathfrak{P}^{(2)}$ has a   second-order echelon on its first step. Then it could be the case that either~(1.1.1) $\mathfrak{P}^{(2)}$ is a   second-order echelon, or~(1.1.2) $\mathfrak{P}^{(2)}$ is a second-order path of length strictly greater than one containing a   second-order echelon on its first step.

If~(1.1.1) then, since $(\mathfrak{P}^{(2)},\mathfrak{Q}^{(2)})\in\mathrm{Ker}(\mathrm{CH}^{(2)})_{s}$ and $\mathrm{CH}^{(2)}_{s}(\mathfrak{P}^{(2)})\in\eta^{(2,\mathcal{A}^{(2)})}[\mathcal{A}^{(2)}]_{s}$, we have, by Lemma~\ref{LDCHDEch}, that $\mathfrak{Q}^{(2)}$ is a   second-order echelon. Thus, the second-order paths $\mathfrak{P}'^{(2)}\circ^{1\mathbf{Pth}_{\boldsymbol{\mathcal{A}}^{(2)}}}_{s}\mathfrak{P}^{(2)}$ and $\mathfrak{Q}'^{(2)}\circ^{1\mathbf{Pth}_{\boldsymbol{\mathcal{A}}^{(2)}}}_{s}\mathfrak{Q}^{(2)}$ are both second-order paths of length strictly greater than one containing a   second-order echelon in its first step, this   second-order echelon being, respectively, $\mathfrak{P}^{(2)}$ and $\mathfrak{Q}^{(2)}$. Then, in virtue of Proposition~\ref{PDCHEchHom}, we have that
\allowdisplaybreaks
\begin{align*}
\mathrm{CH}^{(2)}_{s}\left(
\mathfrak{P}'^{(2)}
\circ^{1\mathbf{Pth}_{\boldsymbol{\mathcal{A}}^{(2)}}}_{s}
\mathfrak{P}^{(2)}
\right)
&=
\mathrm{CH}^{(2)}_{s}\left(
\mathfrak{P}'^{(2)}
\right)
\circ_{s}^{1\mathbf{T}_{\Sigma^{\boldsymbol{\mathcal{A}}^{(2)}}}(X)}
\mathrm{CH}^{(2)}_{s}\left(
\mathfrak{P}^{(2)}
\right)
\\&=
\mathrm{CH}^{(2)}_{s}\left(
\mathfrak{Q}'^{(2)}
\right)
\circ_{s}^{1\mathbf{T}_{\Sigma^{\boldsymbol{\mathcal{A}}^{(2)}}}(X)}
\mathrm{CH}^{(2)}_{s}\left(
\mathfrak{Q}^{(2)}
\right)
\\&=
\mathrm{CH}^{(2)}_{s}\left(
\mathfrak{Q}'^{(2)}
\circ^{1\mathbf{Pth}_{\boldsymbol{\mathcal{A}}^{(2)}}}_{s}
\mathfrak{Q}^{(2)}
\right).
\end{align*}

If~(1.1.2) then the value of the second-order Curry-Howard mapping at $\mathfrak{P}^{(2)}$ is given by
$$
\mathrm{CH}^{(2)}_{s}\left(
\mathfrak{P}^{(2)}
\right)
=
\mathrm{CH}^{(2)}_{s}\left(
\mathfrak{P}^{(2),1,\bb{\mathfrak{P}^{(2)}}-1}
\right)
\circ^{1\mathbf{T}_{\Sigma^{\boldsymbol{\mathcal{A}}^{(2)}}}(X)}_{s}
\mathrm{CH}^{(2)}_{s}\left(
\mathfrak{P}^{(2),0,0}
\right).
$$

Since $(\mathfrak{P}^{(2)},\mathfrak{Q}^{(2)})\in\mathrm{Ker}(\mathrm{CH}^{(2)})_{s}$ and $\mathrm{CH}^{(2)}_{s}(\mathfrak{P}^{(2)})\in\eta^{(2,\mathcal{A}^{(2)})}[\mathcal{A}^{(2)}]_{s}^{\mathrm{int}}$, we have, by Lemma~\ref{LDCHEchInt}, that $\mathfrak{Q}^{(2)}$ is a second-order path of length strictly greater than one containing a   second-order echelon on its first step. The same applies for the second-order path $\mathfrak{Q}'^{(2)}\circ_{s}^{1\mathbf{Pth}_{\boldsymbol{\mathcal{A}}^{(2)}}}\mathfrak{Q}^{(2)}$.

Then, the value of the second-order Curry-Howard mapping at $\mathfrak{Q}^{(2)}$ is given by
$$
\mathrm{CH}^{(2)}_{s}\left(
\mathfrak{Q}^{(2)}
\right)
=
\mathrm{CH}^{(2)}_{s}\left(
\mathfrak{Q}^{(2),1,\bb{\mathfrak{Q}^{(2)}}-1}
\right)
\circ^{1\mathbf{T}_{\Sigma^{\boldsymbol{\mathcal{A}}^{(2)}}}(X)}_{s}
\mathrm{CH}^{(2)}_{s}\left(
\mathfrak{Q}^{(2),0,0}
\right).
$$

Since $(\mathfrak{P}^{(2)},\mathfrak{Q}^{(2)})\in\mathrm{Ker}(\mathrm{CH}^{(2)})_{s}$, we have that 
$
(
\mathfrak{P}^{(2),1,\bb{\mathfrak{P}^{(2)}}-1},
\mathfrak{Q}^{(2),1,\bb{\mathfrak{Q}^{(2)}}-1}
)
$
and 
$(
\mathfrak{P}^{(2),0,0},
\mathfrak{Q}^{(2),0,0}
)
$ 
are pairs in $\mathrm{Ker}(\mathrm{CH}^{(2)})_{s}$. 

Moreover, taking into account that 
\begin{align*}
\mathrm{tg}^{([1],2)}_{s}\left(
\mathfrak{Q}^{(2)}
\right)
&=
\mathrm{tg}^{([1],2)}_{s}\left(
\mathfrak{Q}^{(2),1,\bb{\mathfrak{Q}^{(2)}}-1}
\right),
&
\mathrm{tg}^{([1],2)}_{s}\left(
\mathfrak{P}^{(2)}
\right)
&=
\mathrm{tg}^{([1],2)}_{s}\left(
\mathfrak{P}^{(2),1,\bb{\mathfrak{P}^{(2)}}-1}
\right),
\end{align*}
the following $1$-compositions are well-defined
\begin{align*}
\mathfrak{P}'^{(2)}
\circ^{1\mathbf{Pth}_{\boldsymbol{\mathcal{A}}^{(2)}}}_{s}
\mathfrak{P}^{(2),1,\bb{\mathfrak{P}^{(2)}}-1}
&&
\mathfrak{Q}'^{(2)}
\circ^{1\mathbf{Pth}_{\boldsymbol{\mathcal{A}}^{(2)}}}_{s}
\mathfrak{Q}^{(2),1,\bb{\mathfrak{Q}^{(2)}}-1}.
\end{align*}

Since $(\mathfrak{P}'^{(2)}
\circ^{1\mathbf{Pth}_{\boldsymbol{\mathcal{A}}^{(2)}}}_{s}
\mathfrak{P}^{(2),1,\bb{\mathfrak{P}^{(2)}}-1},s)\prec_{\mathbf{Pth}_{\boldsymbol{\mathcal{A}}^{(2)}}}(\mathfrak{P}'^{(2)}
\circ^{1\mathbf{Pth}_{\boldsymbol{\mathcal{A}}^{(2)}}}_{s}
\mathfrak{P}^{(2)},s)$ and the pairs $
(
\mathfrak{P}^{(2),1,\bb{\mathfrak{P}^{(2)}}-1},
\mathfrak{Q}^{(2),1,\bb{\mathfrak{Q}^{(2)}}-1}
)
$
and 
$(
\mathfrak{P}^{(2),0,0},
\mathfrak{Q}^{(2),0,0}
)
$ 
are in $\mathrm{Ker}(\mathrm{CH}^{(2)})_{s}$, we have, by induction, that 
$$
\left(\mathfrak{P}'^{(2)}
\circ^{1\mathbf{Pth}_{\boldsymbol{\mathcal{A}}^{(2)}}}_{s}
\mathfrak{P}^{(2),1,\bb{\mathfrak{P}^{(2)}}-1},
\mathfrak{Q}'^{(2)}
\circ^{1\mathbf{Pth}_{\boldsymbol{\mathcal{A}}^{(2)}}}_{s}
\mathfrak{Q}^{(2),1,\bb{\mathfrak{Q}^{(2)}}-1}
\right)\in\mathrm{Ker}\left(\mathrm{CH}^{(2)}
\right)_{s}.
$$

So, considering the foregoing, we can affirm that 
\begin{flushleft}
$
\mathrm{CH}^{(2)}_{s}\left(
\mathfrak{P}'^{(2)}
\circ^{1\mathbf{Pth}_{\boldsymbol{\mathcal{A}}^{(2)}}}_{s}
\mathfrak{P}^{(2)}
\right)$
\begin{align*}
\quad
&=
\mathrm{CH}^{(2)}_{s}\left(\left(
\mathfrak{P}'^{(2)}
\circ^{1\mathbf{Pth}_{\boldsymbol{\mathcal{A}}^{(2)}}}_{s}
\mathfrak{P}^{(2)}
\right)^{1,\bb{
\mathfrak{P}'^{(2)}
\circ^{1\mathbf{Pth}_{\boldsymbol{\mathcal{A}}^{(2)}}}_{s}
\mathfrak{P}^{(2)}
}-1}
\right)
\circ^{1\mathbf{T}_{\Sigma^{\boldsymbol{\mathcal{A}}^{(2)}}}(X)}_{s}
\\
&\qquad\qquad\qquad\qquad\qquad\qquad\qquad\qquad
\mathrm{CH}^{(2)}_{s}\left(\left(
\mathfrak{P}'^{(2)}
\circ^{1\mathbf{Pth}_{\boldsymbol{\mathcal{A}}^{(2)}}}_{s}
\mathfrak{P}^{(2)}
\right)^{0,0}
\right)
\\&=
\mathrm{CH}^{(2)}_{s}\left(
\mathfrak{P}'^{(2)}
\circ^{1\mathbf{Pth}_{\boldsymbol{\mathcal{A}}^{(2)}}}_{s}
\mathfrak{P}^{(2),1,\bb{
\mathfrak{P}^{(2)}
}-1}
\right)
\circ^{1\mathbf{T}_{\Sigma^{\boldsymbol{\mathcal{A}}^{(2)}}}(X)}_{s}
\mathrm{CH}^{(2)}_{s}\left(
\mathfrak{P}^{(2),0,0}
\right)
\\&=
\mathrm{CH}^{(2)}_{s}\left(
\mathfrak{Q}'^{(2)}
\circ^{1\mathbf{Pth}_{\boldsymbol{\mathcal{A}}^{(2)}}}_{s}
\mathfrak{Q}^{(2),1,\bb{
\mathfrak{Q}^{(2)}
}-1}
\right)
\circ^{1\mathbf{T}_{\Sigma^{\boldsymbol{\mathcal{A}}^{(2)}}}(X)}_{s}
\mathrm{CH}^{(2)}_{s}\left(
\mathfrak{Q}^{(2),0,0}
\right)
\\&=
\mathrm{CH}^{(2)}_{s}\left(\left(
\mathfrak{Q}'^{(2)}
\circ^{1\mathbf{Pth}_{\boldsymbol{\mathcal{A}}^{(2)}}}_{s}
\mathfrak{Q}^{(2)}
\right)^{1,\bb{
\mathfrak{Q}'^{(2)}
\circ^{1\mathbf{Pth}_{\boldsymbol{\mathcal{A}}^{(2)}}}_{s}
\mathfrak{Q}^{(2)}
}-1}
\right)
\circ^{1\mathbf{T}_{\Sigma^{\boldsymbol{\mathcal{A}}^{(2)}}}(X)}_{s}
\\
&\qquad\qquad\qquad\qquad\qquad\qquad\qquad\qquad
\mathrm{CH}^{(2)}_{s}\left(\left(
\mathfrak{Q}'^{(2)}
\circ^{1\mathbf{Pth}_{\boldsymbol{\mathcal{A}}^{(2)}}}_{s}
\mathfrak{Q}^{(2)}
\right)^{0,0}
\right)
\\&=
\mathrm{CH}^{(2)}_{s}\left(
\mathfrak{Q}'^{(2)}
\circ^{1\mathbf{Pth}_{\boldsymbol{\mathcal{A}}^{(2)}}}_{s}
\mathfrak{Q}^{(2)}
\right).
\end{align*}
\end{flushleft}

Thus case $i=0$ follows.

For case (1.2), i.e., if $i\neq 0$, since $\bb{\mathfrak{P}'^{(2)}
\circ_{s}^{1\mathbf{Pth}_{\boldsymbol{\mathcal{A}}^{(2)}}}
\mathfrak{P}^{(2)}
}=
\bb{\mathfrak{P}'^{(2)}}
+
\bb{\mathfrak{P}^{(2)}}
$, then either (1.2.1) $i\in\bb{\mathfrak{P}^{(2)}}$ or (1.2.2) $i\in [\bb{\mathfrak{P}^{(2)}}, \bb{
\mathfrak{P}'^{(2)}\circ^{1\mathbf{Pth}_{\boldsymbol{\mathcal{A}}^{(2)}}}_{s}\mathfrak{P}^{(2)}
}-1]$.

If (1.2.1), i.e., $i\neq 0$ and $i\in\bb{\mathfrak{P}^{(2)}}$, then $\mathfrak{P}^{(2)}$ is a second-order path of length strictly greater than one containing a   second-order echelon on a step different from the initial one.

Then the value of the second-order Curry-Howard mapping at $\mathfrak{P}^{(2)}$ is given by
$$
\mathrm{CH}^{(2)}_{s}\left(
\mathfrak{P}^{(2)}
\right)
=
\mathrm{CH}^{(2)}_{s}\left(
\mathfrak{P}^{(2),i,\bb{\mathfrak{P}^{(2)}}-1}
\right)
\circ_{s}^{1\mathbf{T}_{\Sigma^{\boldsymbol{\mathcal{A}}^{(2)}}}(X)}
\mathrm{CH}^{(2)}_{s}\left(
\mathfrak{P}^{(2),0,i-1}
\right).
$$

Since $(\mathfrak{P}^{(2)},\mathfrak{Q}^{(2)})\in\mathrm{Ker}(\mathrm{CH}^{(2)})_{s}$ and $\mathrm{CH}^{(2)}_{s}(\mathfrak{P}^{(2)})\in\eta^{(2,\mathcal{A}^{(2)})}[\mathcal{A}^{(2)}]^{\neg\mathrm{int}}_{s}$, we have, by  Lemma~\ref{LDCHEchNInt}, that $\mathfrak{Q}^{(2)}$ is a second-order path of length strictly greater than one containing a   second-order echelon on a step different from the initial one. The same applies for the second-order path $\mathfrak{Q}'^{(2)}\circ_{s}^{1\mathbf{Pth}_{\boldsymbol{\mathcal{A}}^{(2)}}}\mathfrak{Q}^{(2)}$.

Let $j\in\bb{\mathfrak{Q}^{(2)}}$ be the first index for which $\mathfrak{Q}^{(2),j,j}$ is a   second-order echelon. By the previous discussion $j\neq 0$. Then the value of the second-order Curry-Howard mapping at $\mathfrak{Q}^{(2)}$ is given by
$$
\mathrm{CH}^{(2)}_{s}\left(
\mathfrak{Q}^{(2)}
\right)
=
\mathrm{CH}^{(2)}_{s}\left(
\mathfrak{Q}^{(2),j,\bb{\mathfrak{P}^{(2)}}-1}
\right)
\circ_{s}^{1\mathbf{T}_{\Sigma^{\boldsymbol{\mathcal{A}}^{(2)}}}(X)}
\mathrm{CH}^{(2)}_{s}\left(
\mathfrak{Q}^{(2),0,j-1}
\right).
$$

Since $(\mathfrak{P}^{(2)},\mathfrak{Q}^{(2)})\in\mathrm{Ker}(\mathrm{CH}^{(2)})_{s}$, we have that $(\mathfrak{P}^{(2),i,\bb{\mathfrak{P}^{(2)}}-1}, \mathfrak{Q}^{(2),j,\bb{\mathfrak{Q}^{(2)}}-1})$ and $(\mathfrak{P}^{(2),0,i-1},\mathfrak{Q}^{(2),0,j-1})$ are pairs in $\mathrm{Ker}(\mathrm{CH}^{(2)})_{s}$. Moreover, taking into account that
\begin{align*}
\mathrm{tg}^{([1],2)}_{s}\left(\mathfrak{Q}^{(2)}
\right)
&=
\mathrm{tg}^{([1],2)}_{s}\left(
\mathfrak{Q}^{(2),j,\bb{\mathfrak{Q}^{(2)}}-1}
\right),
&
\mathrm{tg}^{([1],2)}_{s}\left(\mathfrak{P}^{(2)}
\right)
&=
\mathrm{tg}^{([1],2)}_{s}\left(
\mathfrak{P}^{(2),i,\bb{\mathfrak{P}^{(2)}}-1}
\right),
\end{align*}
the following $1$-compositions are well-defined
\begin{align*}
\mathfrak{P}'^{(2)}
\circ_{s}^{1\mathbf{Pth}_{\boldsymbol{\mathcal{A}}^{(2)}}}
\mathfrak{P}^{(2),i,\bb{\mathfrak{P}^{(2)}}-1},
&&
\mathfrak{Q}'^{(2)}
\circ_{s}^{1\mathbf{Pth}_{\boldsymbol{\mathcal{A}}^{(2)}}}
\mathfrak{Q}^{(2),j,\bb{\mathfrak{Q}^{(2)}}-1}.
\end{align*}

Since $(\mathfrak{P}'^{(2)}\circ_{s}^{1\mathbf{Pth}_{\boldsymbol{\mathcal{A}}^{(2)}}}
\mathfrak{P}^{(2),i,\bb{\mathfrak{P}^{(2)}}-1},s)
\prec_{\mathbf{Pth}_{\boldsymbol{\mathcal{A}}^{(2)}}}
(\mathfrak{P}'^{(2)}\circ_{s}^{1\mathbf{Pth}_{\boldsymbol{\mathcal{A}}^{(2)}}}
\mathfrak{P}^{(2)},s)
$ and the pairs 
$(
\mathfrak{P}^{(2),i,\bb{\mathfrak{P}^{(2)}}-1},
\mathfrak{Q}^{(2),j,\bb{\mathfrak{Q}^{(2)}}-1}
)
$ and $(\mathfrak{P}'^{(2)},\mathfrak{Q}'^{(2)})$ are in $\mathrm{Ker}(\mathrm{CH}^{(2)})_{s}$, we have, by induction, that
$$
\left(
\mathfrak{P}'^{(2)}
\circ_{s}^{1\mathbf{Pth}_{\boldsymbol{\mathcal{A}}^{(2)}}}
\mathfrak{P}^{(2),i,\bb{\mathfrak{P}^{(2)}}-1},
\mathfrak{Q}'^{(2)}
\circ_{s}^{1\mathbf{Pth}_{\boldsymbol{\mathcal{A}}^{(2)}}}
\mathfrak{Q}^{(2),j,\bb{\mathfrak{Q}^{(2)}}-1}
\right)\in\mathrm{Ker}\left(\mathrm{CH}^{(2)}
\right)_{s}.
$$

So, considering the foregoing, we can affirm that 
\begin{flushleft}
$
\mathrm{CH}^{(2)}_{s}\left(
\mathfrak{P}'^{(2)}
\circ^{1\mathbf{Pth}_{\boldsymbol{\mathcal{A}}^{(2)}}}_{s}
\mathfrak{P}^{(2)}
\right)$
\begin{align*}
\quad
&=
\mathrm{CH}^{(2)}_{s}\left(\left(
\mathfrak{P}'^{(2)}
\circ^{1\mathbf{Pth}_{\boldsymbol{\mathcal{A}}^{(2)}}}_{s}
\mathfrak{P}^{(2)}
\right)^{i,\bb{
\mathfrak{P}'^{(2)}
\circ^{1\mathbf{Pth}_{\boldsymbol{\mathcal{A}}^{(2)}}}_{s}
\mathfrak{P}^{(2)}
}-1}
\right)
\circ^{1\mathbf{T}_{\Sigma^{\boldsymbol{\mathcal{A}}^{(2)}}}(X)}_{s}
\\
&\qquad\qquad\qquad\qquad\qquad\qquad\qquad\qquad
\mathrm{CH}^{(2)}_{s}\left(\left(
\mathfrak{P}'^{(2)}
\circ^{1\mathbf{Pth}_{\boldsymbol{\mathcal{A}}^{(2)}}}_{s}
\mathfrak{P}^{(2)}
\right)^{0,i-1}
\right)
\\&=
\mathrm{CH}^{(2)}_{s}\left(
\mathfrak{P}'^{(2)}
\circ^{1\mathbf{Pth}_{\boldsymbol{\mathcal{A}}^{(2)}}}_{s}
\mathfrak{P}^{(2),i,\bb{
\mathfrak{P}^{(2)}
}-1}
\right)
\circ^{1\mathbf{T}_{\Sigma^{\boldsymbol{\mathcal{A}}^{(2)}}}(X)}_{s}
\mathrm{CH}^{(2)}_{s}\left(
\mathfrak{P}^{(2),0,i-1}
\right)
\\&=
\mathrm{CH}^{(2)}_{s}\left(
\mathfrak{Q}'^{(2)}
\circ^{1\mathbf{Pth}_{\boldsymbol{\mathcal{A}}^{(2)}}}_{s}
\mathfrak{Q}^{(2),j,\bb{
\mathfrak{Q}^{(2)}
}-1}
\right)
\circ^{1\mathbf{T}_{\Sigma^{\boldsymbol{\mathcal{A}}^{(2)}}}(X)}_{s}
\mathrm{CH}^{(2)}_{s}\left(
\mathfrak{Q}^{(2),0,j-1}
\right)
\\&=
\mathrm{CH}^{(2)}_{s}\left(\left(
\mathfrak{Q}'^{(2)}
\circ^{1\mathbf{Pth}_{\boldsymbol{\mathcal{A}}^{(2)}}}_{s}
\mathfrak{Q}^{(2)}
\right)^{j,\bb{
\mathfrak{Q}'^{(2)}
\circ^{1\mathbf{Pth}_{\boldsymbol{\mathcal{A}}^{(2)}}}_{s}
\mathfrak{Q}^{(2)}
}-1}
\right)
\circ^{1\mathbf{T}_{\Sigma^{\boldsymbol{\mathcal{A}}^{(2)}}}(X)}_{s}
\\
&\qquad\qquad\qquad\qquad\qquad\qquad\qquad\qquad
\mathrm{CH}^{(2)}_{s}\left(\left(
\mathfrak{Q}'^{(2)}
\circ^{1\mathbf{Pth}_{\boldsymbol{\mathcal{A}}^{(2)}}}_{s}
\mathfrak{Q}^{(2)}
\right)^{0,j-1}
\right)
\\&=
\mathrm{CH}^{(2)}_{s}\left(
\mathfrak{Q}'^{(2)}
\circ^{1\mathbf{Pth}_{\boldsymbol{\mathcal{A}}^{(2)}}}_{s}
\mathfrak{Q}^{(2)}
\right).
\end{align*}
\end{flushleft}

The case $i\in\bb{\mathfrak{P}^{(2)}}$ follows.

If (1.2.2), i.e., $i\neq 0$ and $i\in [\bb{\mathfrak{P}^{(2)}}, \bb{\mathfrak{P}'\circ^{1\mathbf{Pth}_{\boldsymbol{\mathcal{A}}^{(2)}}}
\mathfrak{P}^{(2)}
}-1]$, then $\mathfrak{P}'^{(2)}$ is a non-$(2,[1])$-identity second-order path containing a   second-order echelon, whilst $\mathfrak{P}^{(2)}$ is an echelonless second-order path.

We will distinguish three cases according to whether (1.2.2.1) $\mathfrak{P}'^{(2)}$ is a   second-order echelon; (1.2.2.2) $\mathfrak{P}'^{(2)}$ is a second-order path of length strictly greater than one containing a   second-order echelon on its first step or (1.2.2.3) $\mathfrak{P}'^{(2)}$ is a second-order path of length strictly greater than one containing a   second-order echelon on a step different from zero. These cases can be proven following a similar argument to those three cases presented above. We leave the details for the interested reader.

If~(2), i.e., if $\mathfrak{P}'^{(2)}\circ_{s}^{1\mathbf{Pth}_{\boldsymbol{\mathcal{A}}^{(2)}}}\mathfrak{P}^{(2)}$ is an echelonless second-order path, it could be the case that~(2.1) $\mathfrak{P}'^{(2)}\circ_{s}^{1\mathbf{Pth}_{\boldsymbol{\mathcal{A}}^{(2)}}}\mathfrak{P}^{(2)}$ is an echelonless second-order path that is not head-constant, or~(2.2) $\mathfrak{P}'^{(2)}\circ_{s}^{1\mathbf{Pth}_{\boldsymbol{\mathcal{A}}^{(2)}}}\mathfrak{P}^{(2)}$ is a head-constant  echelonless second-order path that is not coherent, or~(2.3) $\mathfrak{P}'^{(2)}\circ_{s}^{1\mathbf{Pth}_{\boldsymbol{\mathcal{A}}^{(2)}}}\mathfrak{P}^{(2)}$ is a  coherent head-constant echelonless second-order path.

If~(2.1), let $i\in\bb{\mathfrak{P}'^{(2)}\circ_{s}^{1\mathbf{Pth}_{\boldsymbol{\mathcal{A}}^{(2)}}}\mathfrak{P}^{(2)}}$ be the greatest index for which $(\mathfrak{P}'^{(2)}\circ_{s}^{1\mathbf{Pth}_{\boldsymbol{\mathcal{A}}^{(2)}}}\mathfrak{P}^{(2)})^{0,i}$ is a head-constant second-order path. Since $\bb{\mathfrak{P}'^{(2)}\circ_{s}^{1\mathbf{Pth}_{\boldsymbol{\mathcal{A}}^{(2)}}}\mathfrak{P}^{(2)}}=\bb{\mathfrak{P}'^{(2)}}+\bb{\mathfrak{P}^{(2)}}$, we have that either (2.1.1) $i\in\bb{\mathfrak{P}^{(2)}}-1$, (2.1.2) $i=\bb{\mathfrak{P}^{(2)}}-1$, or (2.1.3) $i\in [\bb{\mathfrak{P}^{(2)}}, \bb{\mathfrak{P}'^{(2)}\circ_{s}^{1\mathbf{Pth}_{\boldsymbol{\mathcal{A}}^{(2)}}}\mathfrak{P}^{(2)}}-1]$.

If~(2.1.1), i.e., if we find ourselves in the case where $\mathfrak{P}'^{(2)}\circ_{s}^{1\mathbf{Pth}_{\boldsymbol{\mathcal{A}}^{(2)}}}\mathfrak{P}^{(2)}$ is an echelonless second-order path that is not head-constant and $i\in\bb{\mathfrak{P}^{(2)}}-1$ is the greatest index for which $(\mathfrak{P}'^{(2)}\circ_{s}^{1\mathbf{Pth}_{\boldsymbol{\mathcal{A}}^{(2)}}}\mathfrak{P}^{(2)})^{0,i}$ is head-constant then, regarding the second-order paths $\mathfrak{P}'^{(2)}$ and $\mathfrak{P}^{(2)}$, we have that
\begin{itemize}
\item[(i)] $\mathfrak{P}^{(2)}$ is an echelonless second-order path that is not head-constant and $i\in\bb{\mathfrak{P}^{(2)}}-1$ is the greatest index for which $\mathfrak{P}^{(2),0,i}$ is head-constant;
\item[(ii)] $\mathfrak{P}'^{(2)}$ is an echelonless second-order path.
\end{itemize}

From (i) and taking into account Definition~\ref{DDCH}, we have that the value of the second-order Curry-Howard mapping at $\mathfrak{P}^{(2)}$ is given by
\begin{align*}
\mathrm{CH}^{(2)}_{s}\left(
\mathfrak{P}^{(2)}
\right)
&=
\mathrm{CH}^{(2)}_{s}\left(
\mathfrak{P}^{(2),i+1,\bb{\mathfrak{P}^{(2)}}-1}
\right)
\circ^{1\mathbf{T}_{\Sigma^{\boldsymbol{\mathcal{A}}^{(2)}}}(X)}_{s}
\mathrm{CH}^{(2)}_{s}\left(
\mathfrak{P}^{(2),0,i}
\right).
\end{align*}

Since $(\mathfrak{P}^{(2)},\mathfrak{Q}^{(2)})\in\mathrm{Ker}(\mathrm{CH}^{(2)})$ and, by Lemma~\ref{LDCHNEchNHd}, $\mathrm{CH}^{(2)}_{s}(\mathfrak{P}^{(2)})\in\mathrm{T}_{\Sigma^{\boldsymbol{\mathcal{A}}^{(2)}}}(X)^{\neg\mathsf{HdC}}_{s}$, we have, again by of Lemma~\ref{LDCHNEchNHd}, that
\begin{itemize}
\item[(ii)] $\mathfrak{Q}^{(2)}$ is an echelonless second-order path that is not head-constant. 
\end{itemize}

Let $j\in\bb{\mathfrak{Q}^{(2)}}-1$ be the greatest index for which $\mathfrak{Q}^{(2),0,j}$ is head-constant. Thus, according to Definition~\ref{DDCH}, we have that the value of the second-order Curry-Howard mapping at $\mathfrak{Q}^{(2)}$ is given by
\begin{align*}
\mathrm{CH}^{(2)}_{s}\left(
\mathfrak{Q}^{(2)}
\right)
&=
\mathrm{CH}^{(2)}_{s}\left(
\mathfrak{Q}^{(2),j+1,\bb{\mathfrak{Q}^{(2)}}-1}
\right)
\circ^{1\mathbf{T}_{\Sigma^{\boldsymbol{\mathcal{A}}^{(2)}}}(X)}_{s}
\mathrm{CH}^{(2)}_{s}\left(
\mathfrak{Q}^{(2),0,j}
\right).
\end{align*} 

From the fact that $(\mathfrak{P}^{(2)},\mathfrak{Q}^{(2)})$ is a pair in $\mathrm{Ker}(\mathrm{CH}^{(2)})_{s}$, we infer that the pairs $(\mathfrak{P}^{(2),i+1,\bb{\mathfrak{P}^{(2)}}-1},\mathfrak{Q}^{(2),j+1,\bb{\mathfrak{Q}^{(2)}}-1})$ and $(\mathfrak{P}^{(2),0,i},\mathfrak{Q}^{(2),0,j})$ are in $\mathrm{Ker}(\mathrm{CH}^{(2)})_{s}$. 

Moreover, taking into account that 
\begin{align*}
\mathrm{tg}^{([1],2)}_{s}\left(\mathfrak{P}^{(2)}
\right)
&=
\mathrm{tg}^{([1],2)}_{s}\left(\mathfrak{P}^{(2),i+1,\bb{\mathfrak{P}^{(2)}}-1}
\right);
\\
\mathrm{tg}^{([1],2)}_{s}\left(\mathfrak{Q}^{(2)}
\right)
&=
\mathrm{tg}^{([1],2)}_{s}\left(\mathfrak{Q}^{(2),j+1,\bb{\mathfrak{Q}^{(2)}}-1}
\right),
\end{align*}
the following $1$-compositions are well-defined
\begin{align*}
\mathfrak{P}'^{(2)}
\circ_{s}^{1\mathbf{Pth}_{\boldsymbol{\mathcal{A}}^{(2)}}}
\mathfrak{P}^{(2),i+1,\bb{\mathfrak{P}^{(2)}}-1},
&&
\mathfrak{Q}'^{(2)}
\circ_{s}^{1\mathbf{Pth}_{\boldsymbol{\mathcal{A}}^{(2)}}}
\mathfrak{Q}^{(2),j+1,\bb{\mathfrak{Q}^{(2)}}-1}.
\end{align*}

Since $(\mathfrak{P}'^{(2)}
\circ_{s}^{1\mathbf{Pth}_{\boldsymbol{\mathcal{A}}^{(2)}}}
\mathfrak{P}^{(2),i+1,\bb{\mathfrak{P}^{(2)}}-1},s)\prec_{\mathbf{Pth}_{\boldsymbol{\mathcal{A}}^{(2)}}}
(\mathfrak{P}'^{(2)}
\circ_{s}^{1\mathbf{Pth}_{\boldsymbol{\mathcal{A}}^{(2)}}}
\mathfrak{P}^{(2)},s)
$ and the pairs $(\mathfrak{P}'^{(2)},\mathfrak{Q}'^{(2)})$ and $(\mathfrak{P}^{(2),i+1,\bb{\mathfrak{P}^{(2)}}-1},\mathfrak{Q}^{(2),j+1,\bb{\mathfrak{Q}^{(2)}}-1})$ are in $\mathrm{Ker}(\mathrm{CH}^{(2)})_{s}$, we have, by induction, that
$$
\left(
\mathfrak{P}'^{(2)}
\circ_{s}^{1\mathbf{Pth}_{\boldsymbol{\mathcal{A}}^{(2)}}}
\mathfrak{P}^{(2),i+1,\bb{\mathfrak{P}^{(2)}}-1},
\mathfrak{Q}'^{(2)}
\circ_{s}^{1\mathbf{Pth}_{\boldsymbol{\mathcal{A}}^{(2)}}}
\mathfrak{Q}^{(2),j+1,\bb{\mathfrak{Q}^{(2)}}-1}
\right)\in\mathrm{Ker}\left(\mathrm{CH}^{(2)}
\right)_{s}.
$$

From (ii) and taking into account Lemma~\ref{LDCHNEch}, we have that
$$
\mathrm{CH}^{(2)}_{s}\left(
\mathfrak{P}'^{(2)}
\right)
\not\in\left(
\eta^{(2,\mathcal{A}^{(2)})}\left[\mathcal{A}^{(2)}
\right]^{\mathrm{pct}}_{s}
\cup
\eta^{(2,1)\sharp}\left[\mathrm{PT}_{\boldsymbol{\mathcal{A}}}
\right]_{s}
\right).
$$

Since $(\mathfrak{P}'^{(2)},\mathfrak{Q}'^{(2)})\in\mathrm{Ker}(\mathrm{CH}^{(2)})_{s}$, we have, again by Lemma~\ref{LDCHNEch}, that 
\begin{itemize}
\item[(iv)] $\mathfrak{Q}'^{(2)}$ is an echelonless second-order path of length at least one.
\end{itemize}

From (iii) and (iv), we infer that $\mathfrak{Q}'^{(2)}\circ_{s}^{1\mathbf{Pth}_{\boldsymbol{\mathcal{A}}^{(2)}}}\mathfrak{Q}^{(2)}$ is an echelonless second-order path of length at least one that is not head-constant. Moreover,  $j\in\bb{\mathfrak{Q}^{(2)}}-1$ is the greatest index for which $(\mathfrak{Q}'^{(2)}\circ_{s}^{1\mathbf{Pth}_{\boldsymbol{\mathcal{A}}^{(2)}}}\mathfrak{Q}^{(2)})^{0,j}$ is head-constant.

So, considering the foregoing, we can affirm that
\begin{flushleft}
$
\mathrm{CH}^{(2)}_{s}\left(
\mathfrak{P}'^{(2)}
\circ_{s}^{1\mathbf{Pth}_{\boldsymbol{\mathcal{A}}^{(2)}}}
\mathfrak{P}^{(2)}
\right)$
\allowdisplaybreaks
\begin{align*}
\quad
&=
\mathrm{CH}^{(2)}_{s}\left(
\left(\mathfrak{P}'^{(2)}
\circ_{s}^{1\mathbf{Pth}_{\boldsymbol{\mathcal{A}}^{(2)}}}
\mathfrak{P}^{(2)}
\right)^{i+1,\bb{
\mathfrak{P}'^{(2)}
\circ_{s}^{1\mathbf{Pth}_{\boldsymbol{\mathcal{A}}^{(2)}}}
\mathfrak{P}^{(2)}
}-1}
\right)
\circ_{s}^{1\mathbf{T}_{\Sigma^{\boldsymbol{\mathcal{A}}^{(2)}}}(X)}
\\&\qquad\qquad\qquad\qquad\qquad\qquad\qquad\qquad
\mathrm{CH}^{(2)}_{s}
\left(\left(\mathfrak{P}'^{(2)}
\circ_{s}^{1\mathbf{Pth}_{\boldsymbol{\mathcal{A}}^{(2)}}}
\mathfrak{P}^{(2)}
\right)^{0,i}
\right)
\\&=
\mathrm{CH}^{(2)}_{s}\left(
\mathfrak{P}'^{(2)}
\circ_{s}^{1\mathbf{Pth}_{\boldsymbol{\mathcal{A}}^{(2)}}}
\mathfrak{P}^{(2),i+1,\bb{
\mathfrak{P}^{(2)}
}-1}
\right)
\circ_{s}^{1\mathbf{T}_{\Sigma^{\boldsymbol{\mathcal{A}}^{(2)}}}(X)}
\mathrm{CH}^{(2)}_{s}
\left(
\mathfrak{P}^{(2),0,i}
\right)
\\&=
\mathrm{CH}^{(2)}_{s}\left(
\mathfrak{Q}'^{(2)}
\circ_{s}^{1\mathbf{Pth}_{\boldsymbol{\mathcal{A}}^{(2)}}}
\mathfrak{Q}^{(2),j+1,\bb{
\mathfrak{Q}^{(2)}
}-1}
\right)
\circ_{s}^{1\mathbf{T}_{\Sigma^{\boldsymbol{\mathcal{A}}^{(2)}}}(X)}
\mathrm{CH}^{(2)}_{s}
\left(
\mathfrak{Q}^{(2),0,j}
\right)
\\&=
\mathrm{CH}^{(2)}_{s}\left(
\left(\mathfrak{Q}'^{(2)}
\circ_{s}^{1\mathbf{Pth}_{\boldsymbol{\mathcal{A}}^{(2)}}}
\mathfrak{Q}^{(2)}
\right)^{j+1,\bb{
\mathfrak{Q}'^{(2)}
\circ_{s}^{1\mathbf{Pth}_{\boldsymbol{\mathcal{A}}^{(2)}}}
\mathfrak{Q}^{(2)}
}-1}
\right)
\circ_{s}^{1\mathbf{T}_{\Sigma^{\boldsymbol{\mathcal{A}}^{(2)}}}(X)}
\\&\qquad\qquad\qquad\qquad\qquad\qquad\qquad\qquad
\mathrm{CH}^{(2)}_{s}
\left(\left(\mathfrak{Q}'^{(2)}
\circ_{s}^{1\mathbf{Pth}_{\boldsymbol{\mathcal{A}}^{(2)}}}
\mathfrak{Q}^{(2)}
\right)^{0,j}
\right)
\\&=
\mathrm{CH}^{(2)}_{s}\left(
\mathfrak{Q}'^{(2)}
\circ_{s}^{1\mathbf{Pth}_{\boldsymbol{\mathcal{A}}^{(2)}}}
\mathfrak{Q}^{(2)}
\right).
\end{align*}
\end{flushleft}

The case $i\in\bb{\mathfrak{P}^{(2)}}-1$ follows.

If~(2.1.2), i.e., if we find ourselves in the case where $\mathfrak{P}'^{(2)}\circ_{s}^{1\mathbf{Pth}_{\boldsymbol{\mathcal{A}}^{(2)}}}\mathfrak{P}^{(2)}$ is an echelonless second-order path that is not head-constant and $i=\bb{\mathfrak{P}^{(2)}}-1$ is the greatest index for which $(\mathfrak{P}'^{(2)}\circ_{s}^{1\mathbf{Pth}_{\boldsymbol{\mathcal{A}}^{(2)}}}\mathfrak{P}^{(2)})^{0,i}$ is head-constant, then, regarding the second-order paths $\mathfrak{P}'^{(2)}$ and $\mathfrak{P}^{(2)}$, we have that
\begin{itemize}
\item[(i)] $\mathfrak{P}^{(2)}$ is a head-constant, echelonless second-order path.
\item[(ii)] $\mathfrak{P}'^{(2)}$ is an echelonless second-order path.
\end{itemize}

Therefore, for a unique word $\mathbf{s}\in S^{\star}-\{\lambda\}$, for a unique operation symbol $\tau\in\Sigma^{\boldsymbol{\mathcal{A}}}_{\mathbf{s},s}$, the family of first-order translations occurring in $\mathfrak{P}^{(2)}$ is a family of first-order translations of type $\tau$.

By Lemma~\ref{LDCHNEchHd}, we have that $\mathrm{CH}^{(2)}_{s}(\mathfrak{P}^{(2)})\in\mathcal{T}(\tau,\mathrm{T}_{\Sigma^{\boldsymbol{\mathcal{A}}^{(2)}}}(X))^{\star}$, which is a subset of $\mathrm{T}_{\Sigma^{\boldsymbol{\mathcal{A}}^{(2)}}}(X)^{\mathsf{HdC}}_{s}$. Since $(\mathfrak{P}^{(2)},\mathfrak{Q}^{(2)})$ is in $\mathrm{Ker}(\mathrm{CH}^{(2)})_{s}$, we have, again by Lemma~\ref{LDCHNEchHd}, that
\begin{itemize}
\item[(iii)] $\mathfrak{Q}^{(2)}$ is a head-constant echelonless second-order path associated to the operation symbol $\tau\in\Sigma_{\mathbf{s},s}$.
\end{itemize}

Also, for a unique word $\mathbf{s}'\in S^{\star}-\{\lambda\}$ and for a unique operation symbol $\tau'\in\Sigma^{\boldsymbol{\mathcal{A}}}_{\mathbf{s}',s}$,  the initial first-order translation occurring in $\mathfrak{P}'^{(2)}$  is a first-order translation of type $\tau'$. Note that $\tau\neq\tau'$, otherwise $i\in\bb{\mathfrak{P}'^{(2)}\circ^{1\mathbf{Pth}_{\boldsymbol{\mathcal{A}}^{(2)}}}
\mathfrak{P}^{(2)}
}$ would not be the greatest index for which $(\mathfrak{P}'^{(2)}\circ^{1\mathbf{Pth}_{\boldsymbol{\mathcal{A}}^{(2)}}}
\mathfrak{P}^{(2)})^{0,i}$ would be head-constant.

By Lemma~\ref{LDCHNEch}, we have that 
$$\mathrm{CH}^{(2)}_{s}\left(\mathfrak{P}'^{(2)}\right)\not
\in
\left(
\eta^{(2,1)\sharp}\left[\mathrm{PT}_{\boldsymbol{\mathcal{A}}}
\right]_{s}
\cup
\eta^{(2,\mathcal{A}^{(2)})}
\left[\mathcal{A}^{(2)}
\right]^{\mathrm{pct}}_{s}
\right).
$$

Since $(\mathfrak{P}'^{(2)},\mathfrak{Q}'^{(2)})$ is a pair in $\mathrm{Ker}(\mathrm{CH}^{(2)})_{s}$, we have, again by Lemma~\ref{LDCHNEch}, that 
\begin{itemize}
\item[(iv)] $\mathfrak{Q}'^{(2)}$ is an echelonless second-order path whose initial first-order translation is associated to the operation symbol $\tau'\in\Sigma^{\boldsymbol{\mathcal{A}}}_{\mathbf{s},s}$.
\end{itemize}

From (iii) and (iv), we infer that $\mathfrak{Q}'^{(2)}\circ_{s}^{1\mathbf{Pth}_{\boldsymbol{\mathcal{A}}^{(2)}}}\mathfrak{Q}^{(2)}$ is an echelonless second-order path that is not head-constant. Moreover,  $j=\bb{\mathfrak{Q}^{(2)}}-1$ is the greatest index for which $(\mathfrak{Q}'^{(2)}\circ_{s}^{1\mathbf{Pth}_{\boldsymbol{\mathcal{A}}^{(2)}}}\mathfrak{Q}^{(2)})^{0,j}$ is head-constant.

So, considering the foregoing, we can affirm that
\begin{flushleft}
$\mathrm{CH}^{(2)}_{s}\left(
\mathfrak{P}'^{(2)}
\circ^{1\mathbf{Pth}_{\boldsymbol{\mathcal{A}}^{(2)}}}_{s}
\mathfrak{P}^{(2)}
\right)
$
\allowdisplaybreaks
\begin{align*}
\quad
&=
\mathrm{CH}^{(2)}_{s}\left(
\left(
\mathfrak{P}'^{(2)}
\circ^{1\mathbf{Pth}_{\boldsymbol{\mathcal{A}}^{(2)}}}_{s}
\mathfrak{P}^{(2)}
\right)^{i+1,\bb{
\mathfrak{P}'^{(2)}
\circ^{1\mathbf{Pth}_{\boldsymbol{\mathcal{A}}^{(2)}}}_{s}
\mathfrak{P}^{(2)}
}-1}
\right)
\circ^{1\mathbf{T}_{\Sigma^{\boldsymbol{\mathcal{A}}^{(2)}}}(X)}_{s}
\\&\qquad\qquad\qquad\qquad\qquad\qquad\qquad\qquad
\mathrm{CH}^{(2)}_{s}\left(
\left(
\mathfrak{P}'^{(2)}
\circ^{1\mathbf{Pth}_{\boldsymbol{\mathcal{A}}^{(2)}}}_{s}
\mathfrak{P}^{(2)}
\right)^{0,i}
\right)
\\&=
\mathrm{CH}^{(2)}_{s}\left(
\mathfrak{P}'^{(2)}
\right)
\circ^{1\mathbf{T}_{\Sigma^{\boldsymbol{\mathcal{A}}^{(2)}}}(X)}_{s}
\mathrm{CH}^{(2)}_{s}\left(
\mathfrak{P}^{(2)}
\right)
\\&=
\mathrm{CH}^{(2)}_{s}\left(
\mathfrak{Q}'^{(2)}
\right)
\circ^{1\mathbf{T}_{\Sigma^{\boldsymbol{\mathcal{A}}^{(2)}}}(X)}_{s}
\mathrm{CH}^{(2)}_{s}\left(
\mathfrak{Q}^{(2)}
\right)
\\
&=
\mathrm{CH}^{(2)}_{s}\left(
\left(
\mathfrak{Q}'^{(2)}
\circ^{1\mathbf{Pth}_{\boldsymbol{\mathcal{A}}^{(2)}}}_{s}
\mathfrak{Q}^{(2)}
\right)^{j+1,\bb{
\mathfrak{Q}'^{(2)}
\circ^{1\mathbf{Pth}_{\boldsymbol{\mathcal{A}}^{(2)}}}_{s}
\mathfrak{Q}^{(2)}
}-1}
\right)
\circ^{1\mathbf{T}_{\Sigma^{\boldsymbol{\mathcal{A}}^{(2)}}}(X)}_{s}
\\&\qquad\qquad\qquad\qquad\qquad\qquad\qquad\qquad
\mathrm{CH}^{(2)}_{s}\left(
\left(
\mathfrak{Q}'^{(2)}
\circ^{1\mathbf{Pth}_{\boldsymbol{\mathcal{A}}^{(2)}}}_{s}
\mathfrak{Q}^{(2)}
\right)^{0,j}
\right)
\\&=
\mathrm{CH}^{(2)}_{s}\left(
\mathfrak{Q}'^{(2)}
\circ^{1\mathbf{Pth}_{\boldsymbol{\mathcal{A}}^{(2)}}}_{s}
\mathfrak{Q}^{(2)}
\right).
\end{align*}
\end{flushleft}

The case $i=\bb{\mathfrak{P}^{(2)}}-1$ follows.

If~(2.1.3), i.e., if we find ourselves in the case where $\mathfrak{P}'^{(2)}\circ_{s}^{1\mathbf{Pth}_{\boldsymbol{\mathcal{A}}^{(2)}}}\mathfrak{P}^{(2)}$ is an echelonless second-order path that is not head-constant and $i\in[\bb{\mathfrak{P}^{(2)}}, \bb{\mathfrak{P}'^{(2)}\circ_{s}^{1\mathbf{Pth}_{\boldsymbol{\mathcal{A}}^{(2)}}}\mathfrak{P}^{(2)}}-1]$ is the greatest index for which $(\mathfrak{P}'^{(2)}\circ_{s}^{1\mathbf{Pth}_{\boldsymbol{\mathcal{A}}^{(2)}}}\mathfrak{P}^{(2)})^{0,i}$ is head-constant then, regarding the second-order paths $\mathfrak{P}'^{(2)}$ and $\mathfrak{P}^{(2)}$, we have that
\begin{itemize}
\item[(i)] $\mathfrak{P}^{(2)}$ is a head-constant echelonless second-order path.
\item[(ii)] $\mathfrak{P}'^{(2)}$ is an echelonless second-order path that is not head-constant and $i-\bb{\mathfrak{P}^{(2)}}\in\bb{\mathfrak{P}'^{(2)}}-1$ is the greatest index for which $\mathfrak{P}'^{(2),0,i-\bb{\mathfrak{P}^{(2)}}}$ is head-constant.
\end{itemize}

Therefore, for a unique word $\mathbf{s}\in S^{\star}-\{\lambda\}$ and a unique operation symbol $\tau\in \Sigma^{\boldsymbol{\mathcal{A}}}_{\mathbf{s},s}$, it is the case that the family of first-order translations occurring in $\mathfrak{P}^{(2)}$ has type $\tau$.

By Lemma~\ref{LDCHNEchHd}, we have that $\mathrm{CH}^{(2)}_{s}(\mathfrak{P}^{(2)})\in\mathcal{T}
(
\tau,\mathrm{T}_{\Sigma^{\boldsymbol{\mathcal{A}}^{(2)}}}(X)
)^{\star}$, which is a subset of $[\mathrm{T}_{\Sigma^{\boldsymbol{\mathcal{A}}^{(2)}}}(X)]^{\mathsf{HdC}}_{s}$. Since $(\mathfrak{P}^{(2)},\mathfrak{Q}^{(2)})$ is in $\mathrm{Ker}(\mathrm{CH}^{(2)})_{s}$, we have, again by Lemma~\ref{LDCHNEchHd}, that
\begin{itemize}
\item[(iii)] $\mathfrak{Q}^{(2)}$ is a head-constant echelonless second-order path associated to the operation symbol $\tau\in\Sigma^{\boldsymbol{\mathcal{A}}}_{\mathbf{s},s}$.
\end{itemize}

Moreover, for the unique word $\mathbf{s}\in S^{\star}-\{\lambda\}$ and the unique operation symbol $\tau\in\Sigma^{\boldsymbol{\mathcal{A}}}_{\mathbf{s},s}$, the initial $i-\bb{\mathfrak{P}^{(2)}}$ first-order translations occurring in $\mathfrak{P}'^{(2)}$ also have  type $\tau$. Note that the operation symbol $\tau$ is the same as in case (i), since $i \in 
[\bb{\mathfrak{P}^{(2)}},\bb{\mathfrak{P}'^{(2)}\circ^{1\mathbf{Pth}_{\boldsymbol{\mathcal{A}}^{(2)}}}
\mathfrak{P}^{(2)}}-1]
$ is the greatest index for which $(\mathfrak{P}'^{(2)}\circ^{1\mathbf{Pth}_{\boldsymbol{\mathcal{A}}^{(2)}}}
\mathfrak{P}^{(2)})^{0,i}$ is head-constant. To simplify the presentation, we let $k=i-\bb{\mathfrak{P}^{(2)}}$.

From (ii) and taking into account Definition~\ref{DDCH}, we have that the value of the second-order Curry-Howard mapping at $\mathfrak{P}'^{(2)}$ is given by
\begin{align*}
\mathrm{CH}^{(2)}_{s}\left(
\mathfrak{P}'^{(2)}
\right)
&=
\mathrm{CH}^{(2)}_{s}\left(
\mathfrak{P}'^{(2),k+1,\bb{\mathfrak{P}'^{(2)}}-1}
\right)
\circ^{1\mathbf{T}_{\Sigma^{\boldsymbol{\mathcal{A}}^{(2)}}}(X)}_{s}
\mathrm{CH}^{(2)}_{s}\left(
\mathfrak{P}'^{(2),0,k}
\right).
\end{align*}

Since $(\mathfrak{P}'^{(2)},\mathfrak{Q}'^{(2)})\in\mathrm{Ker}(\mathrm{CH}^{(2)})$ and, by Lemma~\ref{LDCHNEchNHd}, $\mathrm{CH}^{(2)}_{s}(\mathfrak{P}^{(2)})\in[\mathrm{T}_{\Sigma^{\boldsymbol{\mathcal{A}}^{(2)}}}(X)]^{\neg\mathsf{HdC}}_{s}$, we have, again in virtue of Lemma~\ref{LDCHNEchNHd}, that
\begin{itemize}
\item[(iv)] $\mathfrak{Q'}^{(2)}$ is an echelonless second-order path that is not head-constant. Moreover, its initial first-order translations are associated to the operation symbol $\tau$. 
\end{itemize}

Let $j\in\bb{\mathfrak{Q}'^{(2)}}-1$ be the greatest index for which $\mathfrak{Q}'^{(2),0,j}$ is head-constant. Thus, according to Definition~\ref{DDCH}, we have that the value of the second-order Curry-Howard mapping at $\mathfrak{Q}^{(2)}$ is given by
\begin{align*}
\mathrm{CH}^{(2)}_{s}\left(
\mathfrak{Q}'^{(2)}
\right)
&=
\mathrm{CH}^{(2)}_{s}\left(
\mathfrak{Q}'^{(2),j+1,\bb{\mathfrak{Q}'^{(2)}}-1}
\right)
\circ^{1\mathbf{T}_{\Sigma^{\boldsymbol{\mathcal{A}}^{(2)}}}(X)}_{s}
\mathrm{CH}^{(2)}_{s}\left(
\mathfrak{Q}'^{(2),0,j}
\right).
\end{align*}

From the fact that $(\mathfrak{P}'^{(2)},\mathfrak{Q}'^{(2)})$ is a pair in $\mathrm{Ker}(\mathrm{CH}^{(2)})_{s}$, we infer that the pairs $(\mathfrak{P}'^{(2),k+1,\bb{\mathfrak{P}'^{(2)}}-1},\mathfrak{Q}'^{(2),j+1,\bb{\mathfrak{Q}'^{(2)}}-1})$ and $(\mathfrak{P}'^{(2),0,k},\mathfrak{Q}'^{(2),0,j})$ are in $\mathrm{Ker}(\mathrm{CH}^{(2)})_{s}$. 

Moreover, taking into account that 
\begin{align*}
\mathrm{sc}^{([1],2)}_{s}\left(\mathfrak{P}'^{(2)}\right)
&=
\mathrm{sc}^{([1],2)}_{s}\left(\mathfrak{P}'^{(2),0,k}\right);
&
\mathrm{sc}^{([1],2)}_{s}\left(\mathfrak{Q}'^{(2)}\right)
&=
\mathrm{sc}^{([1],2)}_{s}\left(\mathfrak{Q}'^{(2),0,j}\right),
\end{align*}
the following $1$-compositions are well-defined
\begin{align*}
\mathfrak{P}'^{(2),0,k}
\circ_{s}^{1\mathbf{Pth}_{\boldsymbol{\mathcal{A}}^{(2)}}}
\mathfrak{P}^{(2)},
&&
\mathfrak{Q}'^{(2),0,j}
\circ_{s}^{1\mathbf{Pth}_{\boldsymbol{\mathcal{A}}^{(2)}}}
\mathfrak{Q}^{(2)}.
\end{align*}

Since $(\mathfrak{P}'^{(2),0,k}
\circ_{s}^{1\mathbf{Pth}_{\boldsymbol{\mathcal{A}}^{(2)}}}
\mathfrak{P}^{(2)},s)\prec_{\mathbf{Pth}_{\boldsymbol{\mathcal{A}}^{(2)}}}
(\mathfrak{P}'^{(2)}
\circ_{s}^{1\mathbf{Pth}_{\boldsymbol{\mathcal{A}}^{(2)}}}
\mathfrak{P}^{(2)},s)
$ and the pairs $(\mathfrak{P}^{(2)},\mathfrak{Q}^{(2)})$ and $(\mathfrak{P}'^{(2),0,k},\mathfrak{Q}'^{(2),0,j})$ are in $\mathrm{Ker}(\mathrm{CH}^{(2)})_{s}$, we have, by induction, that
$$
\left(
\mathfrak{P}'^{(2),0,k}
\circ_{s}^{1\mathbf{Pth}_{\boldsymbol{\mathcal{A}}^{(2)}}}
\mathfrak{P}^{(2)},
\mathfrak{Q}'^{(2),0,j}
\circ_{s}^{1\mathbf{Pth}_{\boldsymbol{\mathcal{A}}^{(2)}}}
\mathfrak{Q}^{(2)}
\right)\in\mathrm{Ker}\left(\mathrm{CH}^{(2)}
\right)_{s}.
$$

From (iii) and (iv), we infer that $\mathfrak{Q}'^{(2)}\circ_{s}^{1\mathbf{Pth}_{\boldsymbol{\mathcal{A}}^{(2)}}}\mathfrak{Q}^{(2)}$ is an echelonless second-order path of length at least one that is not head-constant. Moreover,  $l=j+\bb{\mathfrak{Q}^{(2)}}
\in[\bb{\mathfrak{Q}^{(2)}},
\bb{\mathfrak{Q}'^{(2)}\circ_{s}^{1\mathbf{Pth}_{\boldsymbol{\mathcal{A}}^{(2)}}}\mathfrak{Q}^{(2)}}-1
]$ is the greatest index for which $(\mathfrak{Q}'^{(2)}\circ_{s}^{1\mathbf{Pth}_{\boldsymbol{\mathcal{A}}^{(2)}}}\mathfrak{Q}^{(2)})^{0,l}$ is head-constant.

So, considering the foregoing, we can affirm that
\begin{flushleft}
$
\mathrm{CH}^{(2)}_{s}\left(
\mathfrak{P}'^{(2)}
\circ_{s}^{1\mathbf{Pth}_{\boldsymbol{\mathcal{A}}^{(2)}}}
\mathfrak{P}^{(2)}
\right)$
\allowdisplaybreaks
\begin{align*}
\quad
&=
\mathrm{CH}^{(2)}_{s}\left(
\left(\mathfrak{P}'^{(2)}
\circ_{s}^{1\mathbf{Pth}_{\boldsymbol{\mathcal{A}}^{(2)}}}
\mathfrak{P}^{(2)}
\right)^{i+1,\bb{
\mathfrak{P}'^{(2)}
\circ_{s}^{1\mathbf{Pth}_{\boldsymbol{\mathcal{A}}^{(2)}}}
\mathfrak{P}^{(2)}
}-1}
\right)
\circ_{s}^{1\mathbf{T}_{\Sigma^{\boldsymbol{\mathcal{A}}^{(2)}}}(X)}
\\&\qquad\qquad\qquad\qquad\qquad\qquad\qquad\qquad\qquad
\mathrm{CH}^{(2)}_{s}
\left(\left(
\mathfrak{P}'^{(2)}
\circ_{s}^{1\mathbf{Pth}_{\boldsymbol{\mathcal{A}}^{(2)}}}
\mathfrak{P}^{(2)}
\right)^{0,i}
\right)
\\&=
\mathrm{CH}^{(2)}_{s}\left(
\mathfrak{P}'^{(2),k+1,\bb{
\mathfrak{P}'^{(2)}
}-1}
\right)
\circ_{s}^{1\mathbf{T}_{\Sigma^{\boldsymbol{\mathcal{A}}^{(2)}}}(X)}
\mathrm{CH}^{(2)}_{s}
\left(
\mathfrak{P}^{(2),0,k}
\circ_{s}^{1\mathbf{Pth}_{\boldsymbol{\mathcal{A}}^{(2)}}}
\mathfrak{P}^{(2)}
\right)
\\&=
\mathrm{CH}^{(2)}_{s}\left(
\mathfrak{Q}'^{(2),j+1,\bb{
\mathfrak{Q}'^{(2)}
}-1}
\right)
\circ_{s}^{1\mathbf{T}_{\Sigma^{\boldsymbol{\mathcal{A}}^{(2)}}}(X)}
\mathrm{CH}^{(2)}_{s}
\left(
\mathfrak{Q}^{(2),0,j}
\circ_{s}^{1\mathbf{Pth}_{\boldsymbol{\mathcal{A}}^{(2)}}}
\mathfrak{Q}^{(2)}
\right)
\\&=
\mathrm{CH}^{(2)}_{s}\left(
\left(\mathfrak{Q}'^{(2)}
\circ_{s}^{1\mathbf{Pth}_{\boldsymbol{\mathcal{A}}^{(2)}}}
\mathfrak{Q}^{(2)}
\right)^{l+1,\bb{
\mathfrak{Q}'^{(2)}
\circ_{s}^{1\mathbf{Pth}_{\boldsymbol{\mathcal{A}}^{(2)}}}
\mathfrak{Q}^{(2)}
}-1}
\right)
\circ_{s}^{1\mathbf{T}_{\Sigma^{\boldsymbol{\mathcal{A}}^{(2)}}}(X)}
\\&\qquad\qquad\qquad\qquad\qquad\qquad\qquad\qquad\qquad
\mathrm{CH}^{(2)}_{s}
\left(\left(\mathfrak{Q}'^{(2)}
\circ_{s}^{1\mathbf{Pth}_{\boldsymbol{\mathcal{A}}^{(2)}}}
\mathfrak{Q}^{(2)}
\right)^{0,l}
\right)
\\&=
\mathrm{CH}^{(2)}_{s}\left(
\mathfrak{Q}'^{(2)}
\circ_{s}^{1\mathbf{Pth}_{\boldsymbol{\mathcal{A}}^{(2)}}}
\mathfrak{Q}^{(2)}
\right).
\end{align*}
\end{flushleft}

The case $i\in [\bb{\mathfrak{P}^{(2)}},
\bb{
\mathfrak{P}'^{(2)}
\circ_{s}^{1\mathbf{Pth}_{\boldsymbol{\mathcal{A}}^{(2)}}}
\mathfrak{P}^{(2)}
}-1]
$ follows.

This completes case~(2.1).

If~(2.2), i.e., if $\mathfrak{P}'^{(2)}\circ_{s}^{1\mathbf{Pth}_{\boldsymbol{\mathcal{A}}^{(2)}}}\mathfrak{P}^{(2)}$ is a head-constant echelonless second-order path that is not coherent, then let $i\in\bb{\mathfrak{P}'^{(2)}\circ_{s}^{1\mathbf{Pth}_{\boldsymbol{\mathcal{A}}^{(2)}}}\mathfrak{P}^{(2)}}$ be the greatest index for which $(\mathfrak{P}'^{(2)}\circ_{s}^{1\mathbf{Pth}_{\boldsymbol{\mathcal{A}}^{(2)}}}\mathfrak{P}^{(2)})^{0,i}$ is a coherent head-constant echelonless second-order path. Since $\bb{\mathfrak{P}'^{(2)}\circ_{s}^{1\mathbf{Pth}_{\boldsymbol{\mathcal{A}}^{(2)}}}\mathfrak{P}^{(2)}}=\bb{\mathfrak{P}'^{(2)}}+\bb{\mathfrak{P}^{(2)}}$, it follows that either (2.2.1) $i\in\bb{\mathfrak{P}^{(2)}}-1$, (2.2.2) $i=\bb{\mathfrak{P}^{(2)}}-1$, or (2.2.3) $i\in [\bb{\mathfrak{P}^{(2)}}, \bb{\mathfrak{P}'^{(2)}\circ_{s}^{1\mathbf{Pth}_{\boldsymbol{\mathcal{A}}^{(2)}}}\mathfrak{P}^{(2)}}-1]$.

If~(2.2.1), i.e., if we find ourselves in the case where $\mathfrak{P}'^{(2)}\circ_{s}^{1\mathbf{Pth}_{\boldsymbol{\mathcal{A}}^{(2)}}}\mathfrak{P}^{(2)}$ is a head-constant echelonless second-order path  that is not coherent and $i\in\bb{\mathfrak{P}^{(2)}}-1$ is the greatest index for which $(\mathfrak{P}'^{(2)}\circ_{s}^{1\mathbf{Pth}_{\boldsymbol{\mathcal{A}}^{(2)}}}\mathfrak{P}^{(2)})^{0,i}$ is head-constant, then, regarding the second-order paths $\mathfrak{P}'^{(2)}$ and $\mathfrak{P}^{(2)}$, we have that
\begin{itemize}
\item[(i)] $\mathfrak{P}^{(2)}$ is a head-constant echelonless second-order path that is not coherent and $i\in\bb{\mathfrak{P}^{(2)}}-1$ is the greatest index for which $\mathfrak{P}^{(2),0,i}$ is coherent;
\item[(ii)] $\mathfrak{P}'^{(2)}$ is a head-constant echelonless second-order path.
\end{itemize}

Therefore, for a unique word $\mathbf{s}\in S^{\star}-\{\lambda\}$ and  a unique operation symbol $\tau\in\Sigma^{\boldsymbol{\mathcal{A}}}_{\mathbf{s},s}$, the family of first-order translations occurring in $\mathfrak{P}^{(2)}$ is a family of first-order translations of type $\tau$.

From (i) and taking into account Definition~\ref{DDCH}, we have that the value of the second-order Curry-Howard mapping at $\mathfrak{P}^{(2)}$ is given by
\begin{align*}
\mathrm{CH}^{(2)}_{s}\left(
\mathfrak{P}^{(2)}
\right)
&=
\mathrm{CH}^{(2)}_{s}\left(
\mathfrak{P}^{(2),i+1,\bb{\mathfrak{P}^{(2)}}-1}
\right)
\circ^{1\mathbf{T}_{\Sigma^{\boldsymbol{\mathcal{A}}^{(2)}}}(X)}_{s}
\mathrm{CH}^{(2)}_{s}\left(
\mathfrak{P}^{(2),0,i}
\right).
\end{align*}

By Lemma~\ref{LDCHNEchHdNC}, $\mathrm{CH}^{(2)}_{s}(\mathfrak{P}^{(2)})\in \mathcal{T}(\tau,\mathrm{T}_{\Sigma^{\boldsymbol{\mathcal{A}}^{(2)}}}(X))^{+}$, which is a subset of $\mathrm{T}_{\Sigma^{\boldsymbol{\mathcal{A}}^{(2)}}}(X)^{\mathsf{HdC}\And\neg\mathsf{C}}_{s}$. Since $(\mathfrak{P}^{(2)},\mathfrak{Q}^{(2)})\in\mathrm{Ker}(\mathrm{CH}^{(2)})$, we have, again in virtue of Lemma~\ref{LDCHNEchHdNC}, that
\begin{itemize}
\item[(iii)] $\mathfrak{Q}^{(2)}$ is a head-constant echelonless second-order path  that is not coherent associated to the operation symbol $\tau$. 
\end{itemize}

Let $j\in\bb{\mathfrak{Q}^{(2)}}-1$ be the greatest index for which $\mathfrak{Q}^{(2),0,j}$ is head-constant. Thus, according to Definition~\ref{DDCH}, we have that the value of the second-order Curry-Howard mapping at $\mathfrak{Q}^{(2)}$ is given by
\begin{align*}
\mathrm{CH}^{(2)}_{s}\left(
\mathfrak{Q}^{(2)}
\right)
&=
\mathrm{CH}^{(2)}_{s}\left(
\mathfrak{Q}^{(2),j+1,\bb{\mathfrak{Q}^{(2)}}-1}
\right)
\circ^{1\mathbf{T}_{\Sigma^{\boldsymbol{\mathcal{A}}^{(2)}}}(X)}_{s}
\mathrm{CH}^{(2)}_{s}\left(
\mathfrak{Q}^{(2),0,j}
\right).
\end{align*} 

From the fact that $(\mathfrak{P}^{(2)},\mathfrak{Q}^{(2)})$ is a pair in $\mathrm{Ker}(\mathrm{CH}^{(2)})_{s}$, we infer that the pairs $(\mathfrak{P}^{(2),i+1,\bb{\mathfrak{P}^{(2)}}-1},\mathfrak{Q}^{(2),j+1,\bb{\mathfrak{Q}^{(2)}}-1})$ and $(\mathfrak{P}^{(2),0,i},\mathfrak{Q}^{(2),0,j})$ are in $\mathrm{Ker}(\mathrm{CH}^{(2)})_{s}$. 

Moreover, taking into account that 
\begin{align*}
\mathrm{tg}^{([1],2)}_{s}\left(
\mathfrak{P}^{(2)}
\right)
&=
\mathrm{tg}^{([1],2)}_{s}
\left(\mathfrak{P}^{(2),i+1,\bb{\mathfrak{P}^{(2)}}-1}
\right);
&
\mathrm{tg}^{([1],2)}_{s}\left(\mathfrak{Q}^{(2)}
\right)
&=
\mathrm{tg}^{([1],2)}_{s}\left(\mathfrak{Q}^{(2),j+1,\bb{\mathfrak{Q}^{(2)}}-1}
\right),
\end{align*}
the following $1$-compositions are well-defined
\begin{align*}
\mathfrak{P}'^{(2)}
\circ_{s}^{1\mathbf{Pth}_{\boldsymbol{\mathcal{A}}^{(2)}}}
\mathfrak{P}^{(2),i+1,\bb{\mathfrak{P}^{(2)}}-1},
&&
\mathfrak{Q}'^{(2)}
\circ_{s}^{1\mathbf{Pth}_{\boldsymbol{\mathcal{A}}^{(2)}}}
\mathfrak{Q}^{(2),j+1,\bb{\mathfrak{Q}^{(2)}}-1}.
\end{align*}

Since $(\mathfrak{P}'^{(2)}
\circ_{s}^{1\mathbf{Pth}_{\boldsymbol{\mathcal{A}}^{(2)}}}
\mathfrak{P}^{(2),i+1,\bb{\mathfrak{P}^{(2)}}-1},s)\prec_{\mathbf{Pth}_{\boldsymbol{\mathcal{A}}^{(2)}}}
(\mathfrak{P}'^{(2)}
\circ_{s}^{1\mathbf{Pth}_{\boldsymbol{\mathcal{A}}^{(2)}}}
\mathfrak{P}^{(2)},s)
$ and the pairs $(\mathfrak{P}'^{(2)},\mathfrak{Q}'^{(2)})$ and $(\mathfrak{P}^{(2),i+1,\bb{\mathfrak{P}^{(2)}}-1},\mathfrak{Q}^{(2),j+1,\bb{\mathfrak{Q}^{(2)}}-1})$ are in $\mathrm{Ker}(\mathrm{CH}^{(2)})_{s}$, we have, by induction, that
$$
\left(
\mathfrak{P}'^{(2)}
\circ_{s}^{1\mathbf{Pth}_{\boldsymbol{\mathcal{A}}^{(2)}}}
\mathfrak{P}^{(2),i+1,\bb{\mathfrak{P}^{(2)}}-1},
\mathfrak{Q}'^{(2)}
\circ_{s}^{1\mathbf{Pth}_{\boldsymbol{\mathcal{A}}^{(2)}}}
\mathfrak{Q}^{(2),j+1,\bb{\mathfrak{Q}^{(2)}}-1}
\right)\in\mathrm{Ker}
\left(\mathrm{CH}^{(2)}
\right)_{s}.
$$

Moreover, for the unique word $\mathbf{s}\in S^{\star}-\{\lambda\}$ and the unique operation symbol $\tau\in\Sigma^{\boldsymbol{\mathcal{A}}}_{\mathbf{s},s}$, the family of first-order translations occurring in $\mathfrak{P}'^{(2)}$ has type $\tau$. Note that the operation symbol $\tau$ is the same as in case (i), since $\mathfrak{P}'^{(2)}\circ^{1\mathbf{Pth}_{\boldsymbol{\mathcal{A}}^{(2)}}}
\mathfrak{P}^{(2)}$ is head-constant. 

By Lemma~\ref{LDCHNEchHd}, $\mathrm{CH}^{(2)}_{s}(\mathfrak{P}'^{(2)})\in \mathcal{T}(\tau,\mathrm{T}_{\Sigma^{\boldsymbol{\mathcal{A}}^{(2)}}}(X))^{\star}$, which is a subset of $\mathrm{T}_{\Sigma^{\boldsymbol{\mathcal{A}}^{(2)}}}(X)^{\mathsf{HdC}}_{s}$. Since $(\mathfrak{P}^{(2)},\mathfrak{Q}^{(2)})\in\mathrm{Ker}(\mathrm{CH}^{(2)})$, we have, again in virtue of Lemma~\ref{LDCHNEchHd}, that
\begin{itemize}
\item[(iv)] $\mathfrak{Q}^{(2)}$ is a head-constant echelonless second-order path associated to the operation symbol $\tau$. 
\end{itemize}

From (ii) and (iv), we infer that $\mathfrak{Q}'^{(2)}\circ_{s}^{1\mathbf{Pth}_{\boldsymbol{\mathcal{A}}^{(2)}}}\mathfrak{Q}^{(2)}$ is a head-constant echelonless second-order path that is not coherent. Moreover,  $j\in\bb{\mathfrak{Q}^{(2)}}-1$ is the greatest index for which $(\mathfrak{Q}'^{(2)}\circ_{s}^{1\mathbf{Pth}_{\boldsymbol{\mathcal{A}}^{(2)}}}\mathfrak{Q}^{(2)})^{0,j}$ is head-constant.

So, considering the foregoing, we can affirm that
\begin{flushleft}
$
\mathrm{CH}^{(2)}_{s}\left(
\mathfrak{P}'^{(2)}
\circ_{s}^{1\mathbf{Pth}_{\boldsymbol{\mathcal{A}}^{(2)}}}
\mathfrak{P}^{(2)}
\right)$
\allowdisplaybreaks
\begin{align*}
\quad
&=
\mathrm{CH}^{(2)}_{s}\left(
\left(\mathfrak{P}'^{(2)}
\circ_{s}^{1\mathbf{Pth}_{\boldsymbol{\mathcal{A}}^{(2)}}}
\mathfrak{P}^{(2)}
\right)^{i+1,\bb{
\mathfrak{P}'^{(2)}
\circ_{s}^{1\mathbf{Pth}_{\boldsymbol{\mathcal{A}}^{(2)}}}
\mathfrak{P}^{(2)}
}-1}
\right)
\circ_{s}^{1\mathbf{T}_{\Sigma^{\boldsymbol{\mathcal{A}}^{(2)}}}(X)}
\\&\qquad\qquad\qquad\qquad\qquad\qquad\qquad\qquad\qquad
\mathrm{CH}^{(2)}_{s}
\left(\left(\mathfrak{P}'^{(2)}
\circ_{s}^{1\mathbf{Pth}_{\boldsymbol{\mathcal{A}}^{(2)}}}
\mathfrak{P}^{(2)}
\right)^{0,i}
\right)
\\&=
\mathrm{CH}^{(2)}_{s}\left(
\mathfrak{P}'^{(2)}
\circ_{s}^{1\mathbf{Pth}_{\boldsymbol{\mathcal{A}}^{(2)}}}
\mathfrak{P}^{(2),i+1,\bb{
\mathfrak{P}^{(2)}
}-1}
\right)
\circ_{s}^{1\mathbf{T}_{\Sigma^{\boldsymbol{\mathcal{A}}^{(2)}}}(X)}
\mathrm{CH}^{(2)}_{s}
\left(
\mathfrak{P}^{(2),0,i}
\right)
\\&=
\mathrm{CH}^{(2)}_{s}\left(
\mathfrak{Q}'^{(2)}
\circ_{s}^{1\mathbf{Pth}_{\boldsymbol{\mathcal{A}}^{(2)}}}
\mathfrak{Q}^{(2),j+1,\bb{
\mathfrak{Q}^{(2)}
}-1}
\right)
\circ_{s}^{1\mathbf{T}_{\Sigma^{\boldsymbol{\mathcal{A}}^{(2)}}}(X)}
\mathrm{CH}^{(2)}_{s}
\left(
\mathfrak{Q}^{(2),0,j}
\right)
\\&=
\mathrm{CH}^{(2)}_{s}\left(
\left(\mathfrak{Q}'^{(2)}
\circ_{s}^{1\mathbf{Pth}_{\boldsymbol{\mathcal{A}}^{(2)}}}
\mathfrak{Q}^{(2)}
\right)^{j+1,\bb{
\mathfrak{Q}'^{(2)}
\circ_{s}^{1\mathbf{Pth}_{\boldsymbol{\mathcal{A}}^{(2)}}}
\mathfrak{Q}^{(2)}
}-1}
\right)
\circ_{s}^{1\mathbf{T}_{\Sigma^{\boldsymbol{\mathcal{A}}^{(2)}}}(X)}
\\&\qquad\qquad\qquad\qquad\qquad\qquad\qquad\qquad\qquad
\mathrm{CH}^{(2)}_{s}
\left(\left(\mathfrak{Q}'^{(2)}
\circ_{s}^{1\mathbf{Pth}_{\boldsymbol{\mathcal{A}}^{(2)}}}
\mathfrak{Q}^{(2)}
\right)^{0,j}
\right)
\\&=
\mathrm{CH}^{(2)}_{s}\left(
\mathfrak{Q}'^{(2)}
\circ_{s}^{1\mathbf{Pth}_{\boldsymbol{\mathcal{A}}^{(2)}}}
\mathfrak{Q}^{(2)}
\right).
\end{align*}
\end{flushleft}

The case $i\in\bb{\mathfrak{P}^{(2)}}-1$ follows.

If~(2.2.2), i.e., if we find ourselves in the case where $\mathfrak{P}'^{(2)}\circ_{s}^{1\mathbf{Pth}_{\boldsymbol{\mathcal{A}}^{(2)}}}\mathfrak{P}^{(2)}$ is a head-constant echelonless second-order path that is not coherent and $i=\bb{\mathfrak{P}^{(2)}}-1$ is the greatest index for which $(\mathfrak{P}'^{(2)}\circ_{s}^{1\mathbf{Pth}_{\boldsymbol{\mathcal{A}}^{(2)}}}\mathfrak{P}^{(2)})^{0,i}$ is coherent,  then, taking into account Lemma~\ref{LTech}, we infer that $\mathfrak{Q}'^{(2)}\circ_{s}^{1\mathbf{Pth}_{\boldsymbol{\mathcal{A}}^{(2)}}}\mathfrak{Q}^{(2)}$ is a head-constant echelonless second-order path that is not coherent. Moreover, $i=\bb{\mathfrak{Q}^{(2)}}-1$ is the greatest index in $\bb{\mathfrak{Q}'^{(2)}\circ_{s}^{1\mathbf{Pth}_{\boldsymbol{\mathcal{A}}^{(2)}}}\mathfrak{Q}^{(2)}}$ for which $(\mathfrak{Q}'^{(2)}\circ_{s}^{1\mathbf{Pth}_{\boldsymbol{\mathcal{A}}^{(2)}}}\mathfrak{Q}^{(2)})^{0,l}$ is coherent.

So, considering the foregoing, we can affirm that
\begin{flushleft}
$\mathrm{CH}^{(2)}_{s}\left(
\mathfrak{P}'^{(2)}\circ_{s}^{1\mathbf{Pth}_{\boldsymbol{\mathcal{A}}^{(2)}}}\mathfrak{P}^{(2)}
\right)
$
\allowdisplaybreaks
\begin{align*}
\quad&=
\mathrm{CH}^{(2)}_{s}\left(
\left(
\mathfrak{P}'^{(2)}\circ_{s}^{1\mathbf{Pth}_{\boldsymbol{\mathcal{A}}^{(2)}}}\mathfrak{P}^{(2)}
\right)^{\bb{\mathfrak{P}^{(2)}},\bb{\mathfrak{P}'^{(2)}\circ_{s}^{1\mathbf{Pth}_{\boldsymbol{\mathcal{A}}^{(2)}}}\mathfrak{P}^{(2)}}-1}
\right)
\circ_{s}^{1\mathbf{T}_{\Sigma^{\boldsymbol{\mathcal{A}}^{(2)}}}(X)}
\\&\qquad\qquad\qquad\qquad\qquad\qquad\qquad\qquad
\mathrm{CH}^{(2)}_{s}\left(
\left(
\mathfrak{P}'^{(2)}\circ_{s}^{1\mathbf{Pth}_{\boldsymbol{\mathcal{A}}^{(2)}}}\mathfrak{P}^{(2)}
\right)^{0,\bb{\mathfrak{P}^{(2)}}-1}
\right)
\\&=
\mathrm{CH}^{(2)}_{s}\left(
\mathfrak{P}'^{(2)}
\right)
\circ_{s}^{1\mathbf{T}_{\Sigma^{\boldsymbol{\mathcal{A}}^{(2)}}}(X)}
\mathrm{CH}^{(2)}_{s}\left(
\mathfrak{P}^{(2)}
\right)
\\&=
\mathrm{CH}^{(2)}_{s}\left(
\mathfrak{Q}'^{(2)}
\right)
\circ_{s}^{1\mathbf{T}_{\Sigma^{\boldsymbol{\mathcal{A}}^{(2)}}}(X)}
\mathrm{CH}^{(2)}_{s}\left(
\mathfrak{Q}^{(2)}
\right)
\\&=
\mathrm{CH}^{(2)}_{s}\left(
\left(
\mathfrak{Q}'^{(2)}\circ_{s}^{1\mathbf{Pth}_{\boldsymbol{\mathcal{A}}^{(2)}}}\mathfrak{Q}^{(2)}
\right)^{\bb{\mathfrak{Q}^{(2)}},\bb{\mathfrak{Q}'^{(2)}\circ_{s}^{1\mathbf{Pth}_{\boldsymbol{\mathcal{A}}^{(2)}}}\mathfrak{Q}^{(2)}}-1}
\right)
\circ_{s}^{1\mathbf{T}_{\Sigma^{\boldsymbol{\mathcal{A}}^{(2)}}}(X)}
\\&\qquad\qquad\qquad\qquad\qquad\qquad\qquad\qquad
\mathrm{CH}^{(2)}_{s}\left(
\left(
\mathfrak{Q}'^{(2)}\circ_{s}^{1\mathbf{Pth}_{\boldsymbol{\mathcal{A}}^{(2)}}}\mathfrak{Q}^{(2)}
\right)^{0,\bb{\mathfrak{Q}^{(2)}}-1}
\right)
\\&=
\mathrm{CH}^{(2)}_{s}\left(
\mathfrak{Q}'^{(2)}\circ_{s}^{1\mathbf{Pth}_{\boldsymbol{\mathcal{A}}^{(2)}}}\mathfrak{Q}^{(2)}
\right).
\end{align*}
\end{flushleft}

The case $i=\bb{\mathfrak{P}^{(2)}}-1$ follows.

If~(2.2.3), i.e., if we find ourselves in the case where $\mathfrak{P}'^{(2)}\circ_{s}^{1\mathbf{Pth}_{\boldsymbol{\mathcal{A}}^{(2)}}}\mathfrak{P}^{(2)}$ is a head-constant echelonless second-order path  that is not coherent and $i\in[\bb{\mathfrak{P}^{(2)}}, \bb{\mathfrak{P}'^{(2)}\circ_{s}^{1\mathbf{Pth}_{\boldsymbol{\mathcal{A}}^{(2)}}}\mathfrak{P}^{(2)}}-1]$ is the greatest index for which $(\mathfrak{P}'^{(2)}\circ_{s}^{1\mathbf{Pth}_{\boldsymbol{\mathcal{A}}^{(2)}}}\mathfrak{P}^{(2)})^{0,i}$ is coherent,  then, regarding the second-order paths $\mathfrak{P}'^{(2)}$ and $\mathfrak{P}^{(2)}$, we have that
\begin{itemize}
\item[(i)] $\mathfrak{P}^{(2)}$ is a  coherent head-constant and echelonless second-order path.
\item[(ii)] $\mathfrak{P}'^{(2)}$ is a head-constant echelonless second-order path  that is not coherent and $i-\bb{\mathfrak{P}^{(2)}}\in\bb{\mathfrak{P}'^{(2)}}-1$ is the greatest index for which $\mathfrak{P}'^{(2),0,i-\bb{\mathfrak{P}^{(2)}}}$ is coherent.
\end{itemize}

Therefore, for a unique word $\mathbf{s}\in S^{\star}-\{\lambda\}$ and a unique operation symbol $\tau\in \Sigma^{\boldsymbol{\mathcal{A}}}_{\mathbf{s},s}$ the family of first-order  translations occurring in $\mathfrak{P}^{(2)}$ is a family of translations of type $\tau$.

By Lemma~\ref{LDCHNEchHdC}, we have that $\mathrm{CH}^{(2)}_{s}(\mathfrak{P}^{(2)})\in\mathcal{T}(\tau,\mathrm{T}_{\Sigma^{\boldsymbol{\mathcal{A}}^{(2)}}}(X))_{1}$, which is a subset of $[\mathrm{T}_{\Sigma^{\boldsymbol{\mathcal{A}}^{(2)}}}(X)]^{\mathsf{HdC}\And\mathsf{C}}_{s}$. Since $(\mathfrak{P}^{(2)},\mathfrak{Q}^{(2)})$ is in $\mathrm{Ker}(\mathrm{CH}^{(2)})_{s}$, we have, again by Lemma~\ref{LDCHNEchHdC}, that
\begin{itemize}
\item[(iii)] $\mathfrak{Q}^{(2)}$ is a  coherent head-constant echelonless second-order path  associated to the operation symbol $\tau$.
\end{itemize}

Moreover, for the unique word $\mathbf{s}\in S^{\star}-\{\lambda\}$ and the unique operation symbol $\tau\in\Sigma^{\boldsymbol{\mathcal{A}}}_{\mathbf{s},s}$, the family of first-order  translations occurring in $\mathfrak{P}'^{(2)}$ is a family of first-order translations of type $\tau$. Note that the operation symbol $\tau$ is the same as in case (i), since $\mathfrak{P}'^{(2)}\circ^{1\mathbf{Pth}_{\boldsymbol{\mathcal{A}}^{(2)}}}
\mathfrak{P}^{(2)})$ is head-constant. To simplify the presentation, we let $k=i-\bb{\mathfrak{P}^{(2)}}$.

From (ii) and taking into account Definition~\ref{DDCH}, we have that the value of the second-order Curry-Howard mapping at $\mathfrak{P}'^{(2)}$ is given by
\begin{align*}
\mathrm{CH}^{(2)}_{s}\left(
\mathfrak{P}'^{(2)}
\right)
&=
\mathrm{CH}^{(2)}_{s}\left(
\mathfrak{P}'^{(2),k+1,\bb{\mathfrak{P}'^{(2)}}-1}
\right)
\circ^{1\mathbf{T}_{\Sigma^{\boldsymbol{\mathcal{A}}^{(2)}}}(X)}_{s}
\mathrm{CH}^{(2)}_{s}\left(
\mathfrak{P}^{(2),0,k}
\right).
\end{align*}

By Lemma~\ref{LDCHNEchHdNC}, $\mathrm{CH}^{(2)}_{s}(\mathfrak{P}^{(2)})\in\mathcal{T}(\tau, \mathrm{T}_{\Sigma^{\boldsymbol{\mathcal{A}}^{(2)}}}(X))^{+}$, which is a subset of $\mathrm{T}_{\Sigma^{\boldsymbol{\mathcal{A}}^{(2)}}}(X)^{\mathsf{HdC}\And\neg\mathsf{C}}_{s}$. Since $(\mathfrak{P}'^{(2)},\mathfrak{Q}'^{(2)})\in\mathrm{Ker}(\mathrm{CH}^{(2)})$, we have, again in virtue of Lemma~\ref{LDCHNEchHdNC}, that
\begin{itemize}
\item[(iv)] $\mathfrak{Q'}^{(2)}$ is a head-constant echelonless second-order path that is not coherent associated to the operation symbol $\tau$. 
\end{itemize}

Let $j\in\bb{\mathfrak{Q}'^{(2)}}-1$ be the greatest index for which $\mathfrak{Q}'^{(2),0,j}$ is coherent. Thus, according to Definition~\ref{DDCH}, we have that the value of the second-order Curry-Howard mapping at $\mathfrak{Q}'^{(2)}$ is given by
\begin{align*}
\mathrm{CH}^{(2)}_{s}\left(
\mathfrak{Q}'^{(2)}
\right)
&=
\mathrm{CH}^{(2)}_{s}\left(
\mathfrak{Q}'^{(2),j+1,\bb{\mathfrak{Q}'^{(2)}}-1}
\right)
\circ^{1\mathbf{T}_{\Sigma^{\boldsymbol{\mathcal{A}}^{(2)}}}(X)}_{s}
\mathrm{CH}^{(2)}_{s}\left(
\mathfrak{Q}'^{(2),0,j}
\right).
\end{align*} 

From the fact that $(\mathfrak{P}'^{(2)},\mathfrak{Q}'^{(2)})$ is a pair in $\mathrm{Ker}(\mathrm{CH}^{(2)})_{s}$, we have that the pairs $(\mathfrak{P}'^{(2),k+1,\bb{\mathfrak{P}'^{(2)}}-1},\mathfrak{Q}'^{(2),j+1,\bb{\mathfrak{Q}'^{(2)}}-1})$ and $(\mathfrak{P}'^{(2),0,k},\mathfrak{Q}'^{(2),0,j})$ are in $\mathrm{Ker}(\mathrm{CH}^{(2)})_{s}$. 

Moreover, taking into account that 
\begin{align*}
\mathrm{sc}^{([1],2)}_{s}\left(\mathfrak{P}'^{(2)}\right)
&=
\mathrm{sc}^{([1],2)}_{s}\left(\mathfrak{P}'^{(2),0,k}\right),
&
\mathrm{sc}^{([1],2)}_{s}\left(\mathfrak{Q}'^{(2)}\right)
&=
\mathrm{sc}^{([1],2)}_{s}\left(\mathfrak{Q}'^{(2),0,j}\right),
\end{align*}
the following $1$-compositions are well-defined
\begin{align*}
\mathfrak{P}'^{(2),0,k}
\circ_{s}^{1\mathbf{Pth}_{\boldsymbol{\mathcal{A}}^{(2)}}}
\mathfrak{P}^{(2)},
&&
\mathfrak{Q}'^{(2),0,j}
\circ_{s}^{1\mathbf{Pth}_{\boldsymbol{\mathcal{A}}^{(2)}}}
\mathfrak{Q}^{(2)}.
\end{align*}

Since $(\mathfrak{P}'^{(2),0,k}
\circ_{s}^{1\mathbf{Pth}_{\boldsymbol{\mathcal{A}}^{(2)}}}
\mathfrak{P}^{(2)},s)\prec_{\mathbf{Pth}_{\boldsymbol{\mathcal{A}}^{(2)}}}
(\mathfrak{P}'^{(2)}
\circ_{s}^{1\mathbf{Pth}_{\boldsymbol{\mathcal{A}}^{(2)}}}
\mathfrak{P}^{(2)},s)
$ and the pairs $(\mathfrak{P}^{(2)},\mathfrak{Q}^{(2)})$ and $(\mathfrak{P}'^{(2),0,k},\mathfrak{Q}'^{(2),0,j})$ are in $\mathrm{Ker}(\mathrm{CH}^{(2)})_{s}$, we have, by induction, that
$$
\left(
\mathfrak{P}'^{(2),0,k}
\circ_{s}^{1\mathbf{Pth}_{\boldsymbol{\mathcal{A}}^{(2)}}}
\mathfrak{P}^{(2)},
\mathfrak{Q}'^{(2),0,j}
\circ_{s}^{1\mathbf{Pth}_{\boldsymbol{\mathcal{A}}^{(2)}}}
\mathfrak{Q}^{(2)}
\right)\in\mathrm{Ker}(\mathrm{CH}^{(2)})_{s}.
$$

From (iii) and (iv), we infer that $\mathfrak{Q}'^{(2)}\circ_{s}^{1\mathbf{Pth}_{\boldsymbol{\mathcal{A}}^{(2)}}}\mathfrak{Q}^{(2)}$ is a head-constant echelonless second-order path of length at least one that is not coherent. Moreover,  $l=j+\bb{\mathfrak{Q}^{(2)}}
\in[\bb{\mathfrak{Q}^{(2)}},
\bb{\mathfrak{Q}'^{(2)}\circ_{s}^{1\mathbf{Pth}_{\boldsymbol{\mathcal{A}}^{(2)}}}\mathfrak{Q}^{(2)}}-1
]$ is the greatest index for which $(\mathfrak{Q}'^{(2)}\circ_{s}^{1\mathbf{Pth}_{\boldsymbol{\mathcal{A}}^{(2)}}}\mathfrak{Q}^{(2)})^{0,l}$ is coherent.

So, considering the foregoing, we can affirm that
\begin{flushleft}
$
\mathrm{CH}^{(2)}_{s}\left(
\mathfrak{P}'^{(2)}
\circ_{s}^{1\mathbf{Pth}_{\boldsymbol{\mathcal{A}}^{(2)}}}
\mathfrak{P}^{(2)}
\right)$
\allowdisplaybreaks
\begin{align*}
\quad
&=
\mathrm{CH}^{(2)}_{s}\left(
\left(\mathfrak{P}'^{(2)}
\circ_{s}^{1\mathbf{Pth}_{\boldsymbol{\mathcal{A}}^{(2)}}}
\mathfrak{P}^{(2)}
\right)^{i+1,\bb{
\mathfrak{P}'^{(2)}
\circ_{s}^{1\mathbf{Pth}_{\boldsymbol{\mathcal{A}}^{(2)}}}
\mathfrak{P}^{(2)}
}-1}
\right)
\circ_{s}^{1\mathbf{T}_{\Sigma^{\boldsymbol{\mathcal{A}}^{(2)}}}(X)}
\\&\qquad\qquad\qquad\qquad\qquad\qquad\qquad\qquad
\mathrm{CH}^{(2)}_{s}
\left(\left(\mathfrak{P}'^{(2)}
\circ_{s}^{1\mathbf{Pth}_{\boldsymbol{\mathcal{A}}^{(2)}}}
\mathfrak{P}^{(2)}
\right)^{0,i}
\right)
\\&=
\mathrm{CH}^{(2)}_{s}\left(
\mathfrak{P}'^{(2),k+1,\bb{
\mathfrak{P}'^{(2)}
}-1}
\right)
\circ_{s}^{1\mathbf{T}_{\Sigma^{\boldsymbol{\mathcal{A}}^{(2)}}}(X)}
\mathrm{CH}^{(2)}_{s}
\left(
\mathfrak{P}^{(2),0,k}
\circ_{s}^{1\mathbf{Pth}_{\boldsymbol{\mathcal{A}}^{(2)}}}
\mathfrak{P}^{(2)}
\right)
\\&=
\mathrm{CH}^{(2)}_{s}\left(
\mathfrak{Q}'^{(2),j+1,\bb{
\mathfrak{Q}'^{(2)}
}-1}
\right)
\circ_{s}^{1\mathbf{T}_{\Sigma^{\boldsymbol{\mathcal{A}}^{(2)}}}(X)}
\mathrm{CH}^{(2)}_{s}
\left(
\mathfrak{Q}^{(2),0,j}
\circ_{s}^{1\mathbf{Pth}_{\boldsymbol{\mathcal{A}}^{(2)}}}
\mathfrak{Q}^{(2)}
\right)
\\&=
\mathrm{CH}^{(2)}_{s}\left(
\left(\mathfrak{Q}'^{(2)}
\circ_{s}^{1\mathbf{Pth}_{\boldsymbol{\mathcal{A}}^{(2)}}}
\mathfrak{Q}^{(2)}
\right)^{l+1,\bb{
\mathfrak{Q}'^{(2)}
\circ_{s}^{1\mathbf{Pth}_{\boldsymbol{\mathcal{A}}^{(2)}}}
\mathfrak{Q}^{(2)}
}-1}
\right)
\circ_{s}^{1\mathbf{T}_{\Sigma^{\boldsymbol{\mathcal{A}}^{(2)}}}(X)}
\\&\qquad\qquad\qquad\qquad\qquad\qquad\qquad\qquad
\mathrm{CH}^{(2)}_{s}
\left(\left(\mathfrak{Q}'^{(2)}
\circ_{s}^{1\mathbf{Pth}_{\boldsymbol{\mathcal{A}}^{(2)}}}
\mathfrak{Q}^{(2)}
\right)^{0,l}
\right)
\\&=
\mathrm{CH}^{(2)}_{s}\left(
\mathfrak{Q}'^{(2)}
\circ_{s}^{1\mathbf{Pth}_{\boldsymbol{\mathcal{A}}^{(2)}}}
\mathfrak{Q}^{(2)}
\right).
\end{align*}
\end{flushleft}

The case $i\in [\bb{\mathfrak{P}^{(2)}},
\bb{
\mathfrak{P}'^{(2)}
\circ_{s}^{1\mathbf{Pth}_{\boldsymbol{\mathcal{A}}^{(2)}}}
\mathfrak{P}^{(2)}
}-1]
$ follows.

This completes case~(2.2).

If~(2.3), i.e., if we find ourselves in the case where $\mathfrak{P}'^{(2)}\circ_{s}^{1\mathbf{Pth}_{\boldsymbol{\mathcal{A}}^{(2)}}}\mathfrak{P}^{(2)}$ is a coherent head-constant echelonless second-order path then, regarding the second-order paths $\mathfrak{P}'^{(2)}$ and $\mathfrak{P}^{(2)}$, we have that
\begin{itemize}
\item[(i)] $\mathfrak{P}^{(2)}$ is a coherent head-constant echelonless second-order path.
\item[(ii)] $\mathfrak{P}'^{(2)}$ is coherent head-constant echelonless second-order path.
\end{itemize}

Therefore, for a unique word $\mathbf{s}\in S^{\star}-\{\lambda\}$ and a unique operation symbol $\tau\in \Sigma^{\boldsymbol{\mathcal{A}}}_{\mathbf{s},s}$, the family of first-order translations occurring in $\mathfrak{P}^{(2)}$ is a family of first-order translations of type $\tau$.

Let $(\mathfrak{P}^{(2)}_{j})_{j\in\bb{\mathbf{s}}}$ be the family of second-order paths we can extract from $\mathfrak{P}^{(2)}$ in virtue of Lemma~\ref{LDPthExtract}. Then, according to Definition~\ref{DDCH}, we have that the value of the second-order Curry-Howard mapping at $\mathfrak{P}^{(2)}$ is given by
$$
\mathrm{CH}^{(2)}_{s}\left(
\mathfrak{P}^{(2)}
\right)
=
\tau^{\mathbf{T}_{\Sigma^{\boldsymbol{\mathcal{A}}^{(2)}}}(X)}
\left(\left(
\mathrm{CH}^{(2)}_{s_{j}}\left(
\mathfrak{P}^{(2)}_{j}
\right)\right)_{j\in\bb{\mathbf{s}}}
\right).
$$

From Lemma~\ref{LDCHNEchHdC}, we have that $\mathrm{CH}^{(2)}_{s}(\mathfrak{P}^{(2)})\in\mathcal{T}(\tau, \mathrm{T}_{\Sigma^{\boldsymbol{\mathcal{A}}^{(2)}}}(X))_{1}$, which is a subset of $\mathrm{T}_{\Sigma^{\boldsymbol{\mathcal{A}}^{(2)}}}(X)^{\mathsf{HdC}\And\mathsf{C}}_{s}$. Since $(\mathfrak{P}^{(2)},\mathfrak{Q}^{(2)})$ is in $\mathrm{Ker}(\mathrm{CH}^{(2)})_{s}$, we have, again by Lemma~\ref{LDCHNEchHdC}, that
\begin{itemize}
\item[(iii)] $\mathfrak{Q}^{(2)}$ is a coherent head-constant echelonless second-order path  associated to the operation symbol $\tau$.
\end{itemize}

Since (iii), let $(\mathfrak{Q}^{(2)}_{j})_{j\in\bb{\mathbf{s}}}$ be the family of second-order paths we can extract from $\mathfrak{Q}^{(2)}$ in virtue of Lemma~\ref{LDPthExtract}. Then, according to Definition~\ref{DDCH}, we have that the value of the second-order Curry-Howard mapping at $\mathfrak{Q}^{(2)}$ is given by
$$
\mathrm{CH}^{(2)}_{s}\left(
\mathfrak{Q}^{(2)}
\right)
=
\tau^{\mathbf{T}_{\Sigma^{\boldsymbol{\mathcal{A}}^{(2)}}}(X)}
\left(\left(
\mathrm{CH}^{(2)}_{s_{j}}\left(
\mathfrak{Q}^{(2)}_{j}
\right)\right)_{j\in\bb{\mathbf{s}}}
\right).
$$

Since $(\mathfrak{P}^{(2)},\mathfrak{Q}^{(2)})$ is a pair in $\mathrm{Ker}(\mathrm{CH}^{(2)})_{s}$, then we have that, for every $j\in\bb{\mathbf{s}}$, it happens that
$$
\left(\mathfrak{P}^{(2)}_{j},
\mathfrak{Q}^{(2)}_{j}\right)
\in\mathrm{Ker}\left(\mathrm{CH}^{(2)}\right)_{s_{j}}.
$$

Moreover, for the unique word $\mathbf{s}\in S^{\star}-\{\lambda\}$ and the unique operation symbol $\tau\in \Sigma^{\boldsymbol{\mathcal{A}}}_{\mathbf{s},s}$, we have that the family of first-order translations occurring in $\mathfrak{P}'^{(2)}$ is a family of first-order translations of type $\tau$. Note that the operation symbol $\tau$ is the same as in case (i), since $\mathfrak{P}'^{(2)}\circ_{s}^{1\mathbf{Pth}_{\boldsymbol{\mathcal{A}}^{(2)}}}\mathfrak{P}^{(2)}$ is head-constant by hypothesis.

Let $(\mathfrak{P}'^{(2)}_{j})_{j\in\bb{\mathbf{s}}}$ be the family of second-order paths we can extract from $\mathfrak{P}'^{(2)}$ in virtue of Lemma~\ref{LDPthExtract}. Then, according to Definition~\ref{DDCH}, we have that the value of the second-order Curry-Howard mapping at $\mathfrak{P}'^{(2)}$ is given by
$$
\mathrm{CH}^{(2)}_{s}\left(
\mathfrak{P}'^{(2)}
\right)
=
\tau^{\mathbf{T}_{\Sigma^{\boldsymbol{\mathcal{A}}^{(2)}}}(X)}
\left(\left(
\mathrm{CH}^{(2)}_{s_{j}}\left(
\mathfrak{P}'^{(2)}_{j}
\right)\right)_{j\in\bb{\mathbf{s}}}
\right).
$$

From Lemma~\ref{LDCHNEchHdC}, we have that $\mathrm{CH}^{(2)}_{s}(\mathfrak{P}'^{(2)})\in\mathcal{T}(\tau, \mathrm{T}_{\Sigma^{\boldsymbol{\mathcal{A}}^{(2)}}}(X))_{1}$, which is a subset of $[\mathrm{T}_{\Sigma^{\boldsymbol{\mathcal{A}}^{(2)}}}(X)]^{\mathrm{DDADDAHdC}\And\mathsf{C}}_{s}$. Since $(\mathfrak{P}'^{(2)},\mathfrak{Q}'^{(2)})$ is in $\mathrm{Ker}(\mathrm{CH}^{(2)})_{s}$, we have, again by Lemma~\ref{LDCHNEchHdC}, that
\begin{itemize}
\item[(iv)] $\mathfrak{Q}'^{(2)}$ is a coherent head-constant  echelonless second-order path of length at least one associated to the operation symbol $\tau$.
\end{itemize}

Since (ii), let $(\mathfrak{Q}'^{(2)}_{j})_{j\in\bb{\mathbf{s}}}$ be the family of second-order paths we can extract from $\mathfrak{Q}'^{(2)}$ in virtue of Lemma~\ref{LDPthExtract}. Then, according to Definition~\ref{DDCH}, we have that the value of the second-order Curry-Howard mapping at $\mathfrak{Q}'^{(2)}$ is given by
$$
\mathrm{CH}^{(2)}_{s}\left(
\mathfrak{Q}'^{(2)}
\right)
=
\tau^{\mathbf{T}_{\Sigma^{\boldsymbol{\mathcal{A}}^{(2)}}}(X)}
\left(\left(
\mathrm{CH}^{(2)}_{s_{j}}\left(
\mathfrak{Q}'^{(2)}_{j}
\right)\right)_{j\in\bb{\mathbf{s}}}
\right).
$$

Since $(\mathfrak{P}'^{(2)},\mathfrak{Q}'^{(2)})$ is a pair in $\mathrm{Ker}(\mathrm{CH}^{(2)})_{s}$, we have, for every $j\in\bb{\mathbf{s}}$, that
$$
\left(\mathfrak{P}'^{(2)}_{j},
\mathfrak{Q}'^{(2)}_{j}\right)
\in\mathrm{Ker}\left(\mathrm{CH}^{(2)}\right)_{s_{j}}.
$$

From (iii) and (iv), we infer that $\mathfrak{Q}'^{(2)}\circ^{1\mathbf{Pth}_{\boldsymbol{\mathcal{A}}^{(2)}}}\mathfrak{Q}^{(2)}$ is a head-constant echelonless second-order path of length at least one. By Lemma~\ref{LTech}, since $\mathfrak{P}'^{(2)}\circ^{1\mathbf{Pth}_{\boldsymbol{\mathcal{A}}^{(2)}}}\mathfrak{P}^{(2)}$ is coherent and $(\mathfrak{P}^{(2)},\mathfrak{Q}^{(2)})$ and  $(\mathfrak{P}'^{(2)},\mathfrak{Q}'^{(2)})$ are in $\mathrm{Ker}(\mathrm{CH}^{(2)})_{s}$, we have that $\mathfrak{Q}'^{(2)}\circ^{1\mathbf{Pth}_{\boldsymbol{\mathcal{A}}^{(2)}}}\mathfrak{Q}^{(2)}$ is also coherent.

Let us consider $((\mathfrak{P}'^{(2)}\circ_{s}^{1\mathbf{Pth}_{\boldsymbol{\mathcal{A}}^{(2)}}}\mathfrak{P}^{(2)})_{j})_{j\in\bb{\mathbf{s}}}$ and  $((\mathfrak{Q}'^{(2)}\circ_{s}^{1\mathbf{Pth}_{\boldsymbol{\mathcal{A}}^{(2)}}}\mathfrak{Q}^{(2)})_{j})_{j\in\bb{\mathbf{s}}}$, the family of second-order paths we can extract, in virtue of Lemma~\ref{LDPthExtract} from $\mathfrak{P}'^{(2)}\circ_{s}^{1\mathbf{Pth}_{\boldsymbol{\mathcal{A}}^{(2)}}}\mathfrak{P}^{(2)}$ and $\mathfrak{Q}'^{(2)}\circ_{s}^{1\mathbf{Pth}_{\boldsymbol{\mathcal{A}}^{(2)}}}\mathfrak{Q}^{(2)}$, respectively. Let us note that, for every $j\in\bb{\mathbf{s}}$, it is the case that
\begin{align*}
\left(\mathfrak{P}'^{(2)}\circ^{1\mathbf{Pth}_{\boldsymbol{\mathcal{A}}^{(2)}}}_{s}\mathfrak{P}^{(2)}
\right)_{j}
&=
\mathfrak{P}'^{(2)}_{j}
\circ^{1\mathbf{Pth}_{\boldsymbol{\mathcal{A}}^{(2)}}}_{s_{j}}
\mathfrak{P}^{(2)}_{j};
\\
\left(\mathfrak{Q}'^{(2)}\circ^{1\mathbf{Pth}_{\boldsymbol{\mathcal{A}}^{(2)}}}_{s}\mathfrak{Q}^{(2)}
\right)_{j}
&=
\mathfrak{Q}'^{(2)}_{j}
\circ^{1\mathbf{Pth}_{\boldsymbol{\mathcal{A}}^{(2)}}}_{s_{j}}
\mathfrak{Q}^{(2)}_{j}.
\end{align*}

For every $j\in\bb{\mathbf{s}}$, $(\mathfrak{P}'^{(2)}_{j}
\circ^{1\mathbf{Pth}_{\boldsymbol{\mathcal{A}}^{(2)}}}_{s_{j}}
\mathfrak{P}^{(2)}_{j},s_{j})$ $\prec_{\mathbf{Pth}_{\boldsymbol{\mathcal{A}}^{(2)}}}
(\mathfrak{P}'^{(2)}
\circ^{1\mathbf{Pth}_{\boldsymbol{\mathcal{A}}^{(2)}}}_{s}
\mathfrak{P}^{(2)},s)
$ and $(\mathfrak{P}'^{(2)}_{j}, \mathfrak{Q}'^{(2)}_{j})$ and $(\mathfrak{P}^{(2)}_{j}, \mathfrak{Q}^{(2)}_{j})$ are pairs in $\mathrm{Ker}(\mathrm{CH}^{(2)})_{s_{j}}$, hence, by induction, we have that, for every $j\in\bb{\mathbf{s}}$,
$$
\left(\mathfrak{P}'^{(2)}_{j}
\circ^{1\mathbf{Pth}_{\boldsymbol{\mathcal{A}}^{(2)}}}_{s_{j}}
\mathfrak{P}^{(2)}_{j},
\mathfrak{Q}'^{(2)}_{j}
\circ^{1\mathbf{Pth}_{\boldsymbol{\mathcal{A}}^{(2)}}}_{s_{j}}
\mathfrak{Q}^{(2)}_{j}
\right)\in\mathrm{Ker}\left(
\mathrm{CH}^{(2)}
\right)_{s_{j}}.
$$

So, considering the foregoing, we can affirm that
\begin{flushleft}
$\mathrm{CH}^{(2)}_{s}\left(
\mathfrak{P}'^{(2)}
\circ^{1\mathbf{Pth}_{\boldsymbol{\mathcal{A}}^{(2)}}}_{s}
\mathfrak{P}^{(2)}
\right)
$
\allowdisplaybreaks
\begin{align*}
\quad
&=
\tau^{\mathbf{T}_{\Sigma^{\boldsymbol{\mathcal{A}}^{(2)}}}(X)}
\left(\left(
\mathrm{CH}^{(2)}_{s_{j}}\left(
\left(\mathfrak{P}'^{(2)}
\circ^{1\mathbf{Pth}_{\boldsymbol{\mathcal{A}}^{(2)}}}_{s}
\mathfrak{P}^{(2)}
\right)_{j}
\right)\right)_{j\in\bb{\mathbf{s}}}
\right)
\\&=
\tau^{\mathbf{T}_{\Sigma^{\boldsymbol{\mathcal{A}}^{(2)}}}(X)}
\left(\left(
\mathrm{CH}^{(2)}_{s_{j}}\left(
\mathfrak{P}'^{(2)}_{j}
\circ^{1\mathbf{Pth}_{\boldsymbol{\mathcal{A}}^{(2)}}}_{s_{j}}
\mathfrak{P}^{(2)}_{j}
\right)\right)_{j\in\bb{\mathbf{s}}}
\right)
\\&=
\tau^{\mathbf{T}_{\Sigma^{\boldsymbol{\mathcal{A}}^{(2)}}}(X)}
\left(\left(
\mathrm{CH}^{(2)}_{s_{j}}\left(
\mathfrak{Q}'^{(2)}_{j}
\circ^{1\mathbf{Pth}_{\boldsymbol{\mathcal{A}}^{(2)}}}_{s_{j}}
\mathfrak{Q}^{(2)}_{j}
\right)\right)_{j\in\bb{\mathbf{s}}}
\right)
\\&=
\tau^{\mathbf{T}_{\Sigma^{\boldsymbol{\mathcal{A}}^{(2)}}}(X)}
\left(\left(
\mathrm{CH}^{(2)}_{s_{j}}\left(
\left(\mathfrak{Q}'^{(2)}
\circ^{1\mathbf{Pth}_{\boldsymbol{\mathcal{A}}^{(2)}}}_{s}
\mathfrak{Q}^{(2)}
\right)_{j}
\right)\right)_{j\in\bb{\mathbf{s}}}
\right)
\\&=
\mathrm{CH}^{(2)}_{s}\left(
\mathfrak{Q}'^{(2)}
\circ^{1\mathbf{Pth}_{\boldsymbol{\mathcal{A}}^{(2)}}}_{s}
\mathfrak{Q}^{(2)}
\right).
\end{align*}
\end{flushleft}

This completes the case (2.3).

This finishes the proof.
\end{proof}

\chapter{
\texorpdfstring
{On the quotient $[\mathrm{Pth}_{\boldsymbol{\mathcal{A}}^{(2)}}]$}
{On the quotient of second-order paths}
}\label{S2F}

This chapter is devoted to the study of the many-sorted quotient $\mathrm{Pth}_{\boldsymbol{\mathcal{A}}^{(2)}}/{\mathrm{Ker}(\mathrm{CH}^{(2)})}$ which, for simplicity, will be denoted by $[\mathrm{Pth}_{\boldsymbol{\mathcal{A}}^{(2)}}]$. Following this simplification the equivalence class of a second-order path $\mathfrak{P}^{(2)}\in\mathrm{Pth}_{\boldsymbol{\mathcal{A}}^{(2)},s}$, will simply be denoted by $[\mathfrak{P}^{(2)}]_{s}$. Since ${\mathrm{Ker}(\mathrm{CH}^{(2)})}$ is a closed $\Sigma^{\boldsymbol{\mathcal{A}}^{(2)}}$-congruence, we can consider the natural structure of many-sorted partial $\Sigma^{\boldsymbol{\mathcal{A}}^{(2)}}$-congruence on the quotient many-sorted set $[\mathrm{Pth}_{\boldsymbol{\mathcal{A}}^{(2)}}]$, that we will denote by $[\mathbf{Pth}_{\boldsymbol{\mathcal{A}}^{(2)}}]$. We then introduce the operations of $([1],2)$-source, $([1],2)$-target and $([1],2)$-identity second-order path in the quotient, that will be denoted, respectively, by $\mathrm{sc}^{([1],[2])}$, $\mathrm{tg}^{([1],[2])} $ and $\mathrm{ip}^{([2],[1])\sharp}$. As it was shown before, these mappings are also $\Sigma^{\boldsymbol{\mathcal{A}}}$-homomorphisms. We then introduce the operations of $(0,2)$-source, $(0,2)$-target and $(0,2)$-identity second-order path in the quotient, that will be denoted, respectively, by $\mathrm{sc}^{(0,[2])}$, $\mathrm{tg}^{(0,[2])} $ and $\mathrm{ip}^{([2],0)\sharp}$. As it was shown before, these mappings are also $\Sigma$-homomorphisms. We conclude this chapter by investigating how the quotient $[\mathrm{Pth}_{\boldsymbol{\mathcal{A}}^{(2)}}]$ is related to the order $\leq_{\mathbf{Pth}_{\boldsymbol{\mathcal{A}}}^{(2)}}$ that we had defined on the many-sorted set of second-order paths. We first introduce the relation $\leq_{[\mathbf{Pth}_{\boldsymbol{\mathcal{A}}^{(2)}}]}$, which determines that two second-order path classes are related as long as we can find representatives in each class for which the inequality, with respect to $\leq_{\mathbf{Pth}_{\boldsymbol{\mathcal{A}}^{(2)}}}$, of the representatives holds. This definition makes it possible to verify that $\leq_{[\mathbf{Pth}_{\boldsymbol{\mathcal{A}}^{(2)}}]}$ is itself an Artinian partial order in the quotient. In addition, the mappings $\mathrm{pr}^{\mathrm{Ker}(\mathrm{CH}^{(2)})}$, of projection to the class, is order-preserving and inversely order-preserving; the mapping $\mathrm{ip}^{([2],[1])\sharp}$, of $(2,[1])$-identity second-order path mapping in the quotient, is an order embedding; and the monomorphic second-order Curry-Howard mapping, defined on the classes, is still order-preserving.


\begin{restatable}{convention}{CDCHClass}
\label{CDCHClass} 
\index{path!second-order!$[\mathfrak{P}^{(2)}]_{s}$}
\index{path!second-order!$[\mathrm{Pth}_{\boldsymbol{\mathcal{A}}^{(2)}}]$}
In order to simplify the presentation, for a sort $s\in S$ and a second-order path $\mathfrak{P}^{(2)}\in\mathrm{Pth}_{\boldsymbol{\mathcal{A}}^{(2)},s}$, we will let $[\mathfrak{P}^{(2)}]^{}_{s}$ stand for $[\mathfrak{P}^{(2)}]_{\mathrm{Ker}(\mathrm{CH}^{(2)})_{s}}$, the $\mathrm{Ker}(\mathrm{CH}^{(2)})_{s}$-class of $\mathfrak{P}^{(2)}$. Moreover, the $S$-sorted quotient $\mathrm{Pth}_{\boldsymbol{\mathcal{A}}^{(2)}}/{\mathrm{Ker}(\mathrm{CH}^{(2)})}$ will simply be denoted by $[\mathrm{Pth}_{\boldsymbol{\mathcal{A}}^{(2)}}]$.
\end{restatable}

\begin{restatable}{definition}{DDCHQuot}
\index{projection!second-order!$\mathrm{pr}^{\mathrm{Ker}(\mathrm{CH}^{(2)})}$}
\index{Curry-Howard!second-order!$\mathrm{CH}^{(2)\mathrm{m}}$}
\label{DDCHQuot} We consider the (epi,mono)-factorization of $\mathrm{CH}^{(2)}$ given by
\begin{enumerate}
\item the projection $\mathrm{pr}^{\mathrm{Ker}(\mathrm{CH}^{(2)})}$ from $\mathrm{Pth}_{\boldsymbol{\mathcal{A}}^{(2)}}$ to $[\mathrm{Pth}_{\boldsymbol{\mathcal{A}}^{(2)}}]$ that, for every sort $s\in S$, maps a second-order path $\mathfrak{P}^{(2)}$ in $\mathrm{Pth}_{\boldsymbol{\mathcal{A}}^{(2)},s}$ to $[\mathfrak{P}^{(2)}]^{}_{s}$, its equivalence class under the kernel of the second-order Curry-Howard mapping; and
\item the embedding $\mathrm{CH}^{(2)\mathrm{m}}$ from $[\mathrm{Pth}_{\boldsymbol{\mathcal{A}}^{(2)}}]$ to $\mathrm{T}_{\Sigma^{\boldsymbol{\mathcal{A}}^{(2)}}}(X)$ that, for every sort $s\in S$, assigns to an equivalence class $[\mathfrak{P}^{(2)}]^{}_{s}$ in $[\mathrm{Pth}_{\boldsymbol{\mathcal{A}}^{(2)}}]_{s}$ the term $\mathrm{CH}^{(2)}_{s}(\mathfrak{P}^{(2)})$, i.e., the value of the second-order Curry-Howard mapping at any equivalence class representative. We will refer to the mapping $\mathrm{CH}^{(2)\mathrm{m}}$ as the \emph{monomorphic second-order Curry-Howard mapping}.
\end{enumerate}

The reader is advised to consult the diagram appearing in Figure~\ref{FDCHQuot}.
\end{restatable}

\begin{figure}
\begin{center}
\begin{tikzpicture}
[ACliment/.style={-{To [angle'=45, length=5.75pt, width=4pt, round]}},scale=.8]
\node[] (2p) at (0,0) [] {$\mathrm{Pth}_{\boldsymbol{\mathcal{A}}^{(2)}}$};
\node[] (2pq) at (6,0) [] {$[\mathrm{Pth}_{\boldsymbol{\mathcal{A}}^{(2)}}]
$};
\node[] (2t) at (0,-3) [] {$\mathrm{T}_{\Sigma^{\boldsymbol{\mathcal{A}}^{(2)}}}(X)$};

\draw[ACliment]  (2p) 	to node [above]	
{$\mathrm{pr}^{\mathrm{Ker}(\mathrm{CH}^{(2)})}$} (2pq);
\draw[ACliment]  (2p) 	to node [left]	
{$\mathrm{CH}^{(2)}$} (2t);
\draw[ACliment, bend left=10]  (2pq) 	to node [below right]	
{$\mathrm{CH}^{(2)\mathrm{m}}$} (2t);

\end{tikzpicture}
\end{center}
\caption{Many-sorted quotient mappings at layer 2.}
\label{FDCHQuot}
\end{figure}

We next introduce the notation for the composition of the embeddings to $\mathrm{Pth}_{\boldsymbol{\mathcal{A}}^{(2)}}$ with the projection to the kernel of the second-order Curry-Howard mapping.

\begin{restatable}{definition}{DDCHEch}
\label{DDCHEch} We will denote by
\begin{enumerate}
\item $\mathrm{ip}^{([2],X)}$ the $S$-sorted mapping from $X$ to $[\mathrm{Pth}_{\boldsymbol{\mathcal{A}}^{(2)}}]$ given by the composition $\mathrm{ip}^{([2]^{},X)}=\mathrm{pr}^{\mathrm{Ker}(\mathrm{CH}^{(2)})}\circ\mathrm{ip}^{(2,X)}$, i.e., for every sort $s\in S$, $\mathrm{ip}^{([2]^{},X)}_{s}$ sends a variable $x\in X_{s}$ to the class $[\mathrm{ip}^{(2,X)}_{s}(x)]^{}_{s}$ in $[\mathrm{Pth}_{\boldsymbol{\mathcal{A}}^{(2)}}]_{s}$.
\index{identity!second-order!$\mathrm{ip}^{([2],X)}$}
\item $\mathrm{ech}^{([2],\mathcal{A})}$ the $S$-sorted mapping from $\mathcal{A}$ to $[\mathrm{Pth}_{\boldsymbol{\mathcal{A}}^{(2)}}]$ given by the composition $\mathrm{ech}^{([2]^{},\mathcal{A})}=\mathrm{pr}^{\mathrm{Ker}(\mathrm{CH}^{(2)})}\circ\mathrm{ech}^{(2,\mathcal{A})}$, i.e., for every sort $s\in S$, $\mathrm{ech}^{([2]^{},\mathcal{A})}_{s}$ sends a rewrite rule $\mathfrak{p}\in \mathcal{A}_{s}$ to the class $[\mathrm{ech}^{(2,\mathcal{A})}_{s}(\mathfrak{p})]^{}_{s}$ in $[\mathrm{Pth}_{\boldsymbol{\mathcal{A}}^{(2)}}]_{s}$.
\index{echelon!second-order!$\mathrm{ech}^{([2],\mathcal{A})}$}
\item $\mathrm{ech}^{([2],\mathcal{A}^{(2)})}$ the $S$-sorted mapping from $\mathcal{A}^{(2)}$ to $[\mathrm{Pth}_{\boldsymbol{\mathcal{A}}^{(2)}}]$ given by the composition $\mathrm{ech}^{([2]^{},\mathcal{A}^{(2)})}=\mathrm{pr}^{\mathrm{Ker}(\mathrm{CH}^{(2)})}\circ\mathrm{ech}^{(2,\mathcal{A}^{(2)})}$, i.e., for every sort $s\in S$, $\mathrm{ech}^{([2]^{},\mathcal{A}^{(2)})}_{s}$ sends a second-order rewrite rule $\mathfrak{p}^{(2)}\in \mathcal{A}^{(2)}_{s}$ to the class $[\mathrm{ech}^{(2,\mathcal{A}^{(2)})}_{s}(\mathfrak{p}^{(2)})]^{}_{s}$ in $[\mathrm{Pth}_{\boldsymbol{\mathcal{A}}^{(2)}}]_{s}$.
\index{echelon!second-order!$\mathrm{ech}^{([2],\mathcal{A}^{(2)})}$}
\end{enumerate}
The above $S$-sorted mappings are depicted in the diagram of Figure~\ref{FDCHEch}.
\end{restatable}

\begin{figure}
\begin{tikzpicture}
[ACliment/.style={-{To [angle'=45, length=5.75pt, width=4pt, round]}},scale=.8]
\node[] (x) at (0,0) [] {$X$};
\node[] (a) at (0,-1.5) [] {$\mathcal{A}$};
\node[] (a2) at (0,-3) [] {$\mathcal{A}^{(2)}$};
\node[] (T) at (6,-3) [] {$
[\mathrm{Pth}_{\boldsymbol{\mathcal{A}}^{(2)}}]$};
\draw[ACliment, bend left=20]  (x) to node [above right] {$\mathrm{ip}^{([2]^{},X)}$} (T);
\draw[ACliment, bend left=10]  (a) to node [midway, fill=white] {$\mathrm{ech}^{([2]^{},\mathcal{A})}$} (T);
\draw[ACliment]  (a2) to node [below] {$\mathrm{ech}^{([2]^{},\mathcal{A}^{(2)})}$} (T);
\end{tikzpicture}
\caption{Quotient path embeddings relative to $X$, $\mathcal{A}$ and $\mathcal{A}^{(2)}$.}\label{FDCHEch}
\end{figure}

Finally, we introduce some terminology for the $S$-sorted mappings of sources, targets and identity second-order paths at each layer with respect to the many-sorted first-order quotient set of second-order paths with respect to the kernel of the second-order Curry-Howard mapping.

\begin{restatable}{definition}{DDCHDU}
\label{DDCHDU} We will denote by
\begin{enumerate}
\item  $\mathrm{sc}^{([1],[2]^{})}$ the $S$-sorted mapping from $[\mathrm{Pth}_{\boldsymbol{\mathcal{A}}^{(2)}}]$ to $[\mathrm{PT}_{\boldsymbol{\mathcal{A}}}]$ that, for every $s\in S$, assigns to an equivalence class $[\mathfrak{P}^{(2)}]^{}_{s}$ in $[\mathrm{Pth}_{\boldsymbol{\mathcal{A}}^{(2)}}]_{s}$ with $\mathfrak{P}^{(2)}\in\mathrm{Pth}_{\boldsymbol{\mathcal{A}}^{(2)},s}$ the path term class $\mathrm{sc}^{([1],2)}_{s}(\mathfrak{P}^{(2)})$, i.e., the value of the mapping $\mathrm{sc}^{([1],2)}_{s}$ at any second-order representative. That $\mathrm{sc}^{([1],[2]^{})}$ is well-defined follows from Lemma~\ref{LDCH};
\index{source!second-order!$\mathrm{sc}^{([1],[2]^{})}$}
\item $\mathrm{tg}^{([1],[2]^{})}$ the $S$-sorted mapping  from $[\mathrm{Pth}_{\boldsymbol{\mathcal{A}}^{(2)}}]$ to $[\mathrm{PT}_{\boldsymbol{\mathcal{A}}}]$ that, for every $s\in S$, assigns to an equivalence class $[\mathfrak{P}^{(2)}]^{}_{s}$ in $[\mathrm{Pth}_{\boldsymbol{\mathcal{A}}^{(2)}}]_{s}$ with $\mathfrak{P}^{(2)}\in\mathrm{Pth}_{\boldsymbol{\mathcal{A}}^{(2)},s}$ the path term class $\mathrm{tg}^{([1],2)}_{s}(\mathfrak{P}^{(2)})$, i.e., the value of the mapping $\mathrm{tg}^{([1],2)}_{s}$ at any second-order representative. That $\mathrm{tg}^{([1],[2]^{})}$ is well-defined follows from Lemma~\ref{LDCH};
\index{target!second-order!$\mathrm{tg}^{([1],[2]^{})}$}
\item $\mathrm{ip}^{([2]^{},[1])\sharp}$  the $S$-sorted mapping from $[\mathrm{PT}_{\boldsymbol{\mathcal{A}}}]$ to $[\mathrm{Pth}_{\boldsymbol{\mathcal{A}}^{(2)}}]$ given by the composition
$\mathrm{ip}^{([2]^{},[1])\sharp}=\mathrm{pr}^{\mathrm{Ker}(\mathrm{CH}^{(2)})}\circ\mathrm{ip}^{(2,[1])\sharp}$. That is, for every sort $s\in S$, $\mathrm{ip}^{([2]^{},[1])\sharp}$ assigns to a path term class $[P]_{s}$ in $[\mathrm{PT}_{\boldsymbol{\mathcal{A}}}]_{s}$ the class with respect to the kernel of the second-order Curry-Howard mapping of the $(2,[1])$-identity second-order path $\mathrm{ip}^{(2,[1])\sharp}_{s}([P]_{s})$, i.e., 
$
\mathrm{ip}^{([2]^{},[1])\sharp}_{s}(
[P]_{s}
)=[\mathrm{ip}^{(2,[1])\sharp}_{s}([P]_{s})]^{}_{s}.
$
\index{identity!second-order!$\mathrm{ip}^{([2]^{},[1])\sharp}$}
\end{enumerate}
\end{restatable}

\begin{proposition} The following equalities hold
\begin{itemize}
\item[(i)] $\mathrm{sc}^{([1],[2]^{})}\circ \mathrm{ip}^{([2]^{},[1])\sharp}=\mathrm{id}^{
[
\mathrm{PT}_{\boldsymbol{\mathcal{A}}}
]
}$;
\item[(ii)] $\mathrm{tg}^{([1],[2]^{})}\circ \mathrm{ip}^{([2]^{},[1])\sharp}=\mathrm{id}^{
[
\mathrm{PT}_{\boldsymbol{\mathcal{A}}}
]
}$.
\end{itemize}
\end{proposition}

\begin{proposition} The following equalities hold
\begin{itemize}
\item[(i)] $\mathrm{sc}^{([1],[2]^{})}\circ \mathrm{ip}^{([2]^{},X)}
=\eta^{([1],X)};$
\item[(ii)] $\mathrm{tg}^{([1],[2]^{})}\circ \mathrm{ip}^{([2]^{},X)}
=\eta^{([1],X)};$
\item[(iii)] $\mathrm{ip}^{([2]^{},[1])\sharp}\circ \eta^{([1],X)}=
\mathrm{ip}^{([2]^{},X)}.
$
\end{itemize}

The reader is advised to consult the diagram appearing in Figure~\ref{FDCHDU}.
\end{proposition}

\begin{figure}
\begin{center}
\begin{tikzpicture}
[ACliment/.style={-{To [angle'=45, length=5.75pt, width=4pt, round]}},scale=1]
\node[] (x) at (0,0) [] {$X$};
\node[] (pt) at (6,0) [] {$[\mathrm{PT}_{\boldsymbol{\mathcal{A}}}]$};
\node[] (2pq) at (6,-3) []  {$[\mathrm{Pth}_{\boldsymbol{\mathcal{A}}^{(2)}}]
$};

\draw[ACliment]  (x) 	to node [above right]	
{$\eta^{([1],X)}$} (pt);
\draw[ACliment, bend right=10]  (x) 	to node [below left]	
{$\mathrm{ip}^{([2]^{},X)}$} (2pq);

\node[] (B1) at (6,-1.5)  [] {};
\draw[ACliment]  ($(B1)+(0,1.2)$) to node [above, fill=white] {
$\mathrm{ip}^{([2]^{},[1])\sharp}$
} ($(B1)+(0,-1.2)$);
\draw[ACliment, bend right]  ($(B1)+(.5,-1.2)$) to node [ below, fill=white] {
$ \mathrm{tg}^{([1],[2]^{})}$
} ($(B1)+(.5,1.2)$);
\draw[ACliment, bend left]  ($(B1)+(-.5,-1.2)$) to node [below, fill=white] {
$ \mathrm{sc}^{([1],[2]^{})}$
} ($(B1)+(-.5,1.2)$);

\end{tikzpicture}
\end{center}
\caption{Quotient mappings relative to $X$ at layers 1 \& 2.}
\label{FDCHDU}
\end{figure}

\begin{proposition} The following equalities hold
\begin{itemize}
\item[(i)] $\mathrm{sc}^{([1],[2]^{})}\circ \mathrm{ech}^{([2]^{},\mathcal{A})}
=\eta^{([1],\mathcal{A})};$
\item[(ii)] $\mathrm{tg}^{([1],[2]^{})}\circ \mathrm{ech}^{([2]^{},\mathcal{A})}
= \eta^{([1],\mathcal{A})};$
\item[(iii)] $\mathrm{ip}^{([2]^{},[1])\sharp}\circ \eta^{([1],\mathcal{A})}=
\mathrm{ech}^{([2]^{},\mathcal{A})}.
$
\end{itemize}

The reader is advised to consult the diagram appearing in Figure~\ref{FDCHDAQ}.
\end{proposition}

\begin{figure}
\begin{center}
\begin{tikzpicture}
[ACliment/.style={-{To [angle'=45, length=5.75pt, width=4pt, round]}},scale=1]
\node[] (x) at (0,0) [] {$\mathcal{A}$};
\node[] (pt) at (6,0) [] {$[\mathrm{PT}_{\boldsymbol{\mathcal{A}}}]$};
\node[] (2pq) at (6,-3) []  {$[\mathrm{Pth}_{\boldsymbol{\mathcal{A}}^{(2)}}]
$};

\draw[ACliment]  (x) 	to node [above right]	
{$\eta^{([1],\mathcal{A})}$} (pt);
\draw[ACliment, bend right=10]  (x) 	to node [below left]	
{$\mathrm{ech}^{([2]^{},\mathcal{A})}$} (2pq);

\node[] (B1) at (6,-1.5)  [] {};
\draw[ACliment]  ($(B1)+(0,1.2)$) to node [above, fill=white] {
$\mathrm{ip}^{([2]^{},[1])\sharp}$
} ($(B1)+(0,-1.2)$);
\draw[ACliment, bend right]  ($(B1)+(.5,-1.2)$) to node [ below, fill=white] {
$ \mathrm{tg}^{([1],[2]^{})}$
} ($(B1)+(.5,1.2)$);
\draw[ACliment, bend left]  ($(B1)+(-.5,-1.2)$) to node [below, fill=white] {
$ \mathrm{sc}^{([1],[2]^{})}$
} ($(B1)+(-.5,1.2)$);

\end{tikzpicture}
\end{center}
\caption{Quotient mappings relative to $\mathcal{A}$ at layers 1 \& 2.}
\label{FDCHDAQ}
\end{figure}

\begin{restatable}{definition}{DDCHDZ}
\label{DDCHDZ} We will denote by
\begin{enumerate}
\item $\mathrm{sc}^{(0,[2]^{})}$ the $S$-sorted mapping from $[\mathrm{Pth}_{\boldsymbol{\mathcal{A}}^{(2)}}]$ to $\mathrm{T}_{\Sigma}(X)$ that, for every $s\in S$, assigns to an equivalence class $[\mathfrak{P}^{(2)}]^{}_{s}$ in $[\mathrm{Pth}_{\boldsymbol{\mathcal{A}}^{(2)}}]_{s}$ with $\mathfrak{P}^{(2)}\in\mathrm{Pth}_{\boldsymbol{\mathcal{A}}^{(2)},s}$ the  term $\mathrm{sc}^{(0,2)}_{s}(\mathfrak{P}^{(2)})$, i.e., the value of the mapping $\mathrm{sc}^{(0,2)}_{s}$ at any second-order representative. That $\mathrm{sc}^{(0,2)}$ is well-defined follows from Corollary~\ref{CDCH};
\index{source!second-order!$\mathrm{sc}^{(0,[2]^{})}$}
\item $\mathrm{tg}^{(0,[2]^{})}$ the $S$-sorted mapping from $[\mathrm{Pth}_{\boldsymbol{\mathcal{A}}^{(2)}}]$ to $\mathrm{T}_{\Sigma}(X)$ that, for every $s\in S$, assigns to an equivalence class $[\mathfrak{P}^{(2)}]^{}_{s}$ in $[\mathrm{Pth}_{\boldsymbol{\mathcal{A}}^{(2)}}]_{s}$ with $\mathfrak{P}^{(2)}\in\mathrm{Pth}_{\boldsymbol{\mathcal{A}}^{(2)},s}$ the  term $\mathrm{tg}^{(0,2)}_{s}(\mathfrak{P}^{(2)})$, i.e., the value of the mapping $\mathrm{tg}^{(0,2)}_{s}$ at any second-order representative. That $\mathrm{tg}^{(0,2)}$ is well-defined follows from Corollary~\ref{CDCH};
\index{target!second-order!$\mathrm{tg}^{(0,[2]^{})}$}
\item $\mathrm{ip}^{([2]^{},0)\sharp}$  the $S$-sorted mapping from $\mathrm{T}_{\Sigma}(X)$ to $[\mathrm{Pth}_{\boldsymbol{\mathcal{A}}^{(2)}}]$ given by the composition
$\mathrm{ip}^{([2]^{},0)\sharp}=\mathrm{pr}^{\mathrm{Ker}(\mathrm{CH}^{(2)})}\circ\mathrm{ip}^{(2,0)\sharp}$. That is, for every sort $s\in S$, $\mathrm{ip}^{([2]^{},0)\sharp}$ assigns to a term $P$ in $\mathrm{T}_{\Sigma}(X)_{s}$ the class with respect to the kernel of the second-order Curry-Howard mapping of the $(2,0)$-identity second-order path $\mathrm{ip}^{(2,0)\sharp}_{s}(P)$, i.e., 
$
\mathrm{ip}^{([2]^{},0)\sharp}_{s}(
P
)=[\mathrm{ip}^{(2,0)\sharp}_{s}(P)]^{}_{s}.
$
\index{identity!second-order!$\mathrm{ip}^{([2]^{},0)\sharp}$}
\end{enumerate} 
\end{restatable}

\begin{proposition} The following equalities hold
\begin{itemize}
\item[(i)] $\mathrm{sc}^{(0,[2]^{})}=\mathrm{sc}^{(0,[1])}\circ \mathrm{sc}^{([1],[2]^{})}$;
\item[(ii)] $\mathrm{tg}^{(0,[2]^{})}=\mathrm{tg}^{(0,[1])}\circ\mathrm{tg}^{([1],[2]^{})}$;
\item[(iii)] $\mathrm{ip}^{([2]^{},0)\sharp}=\mathrm{ip}^{([2]^{},[1])\sharp}\circ\mathrm{ip}^{([1],0)\sharp}$.
\end{itemize}
\end{proposition}

\begin{proposition} The following equalities hold
\begin{itemize}
\item[(i)] $\mathrm{sc}^{(0,[2]^{})}\circ \mathrm{ip}^{([2]^{},0)\sharp}=\mathrm{id}^{
\mathrm{T}_{\Sigma}(X)
}$;
\item[(ii)] $\mathrm{tg}^{(0,[2]^{})}\circ\mathrm{ip}^{([2]^{},0)\sharp}=\mathrm{id}^{
\mathrm{T}_{\Sigma}(X)
}$.
\end{itemize}
\end{proposition}

\begin{proposition} The following equalities hold
\begin{itemize}
\item[(i)] $\mathrm{sc}^{(0,[2]^{})}\circ \mathrm{ip}^{([2]^{},X)}
=\eta^{(0,X)};$
\item[(ii)] $\mathrm{tg}^{(0,[2]^{})}\circ \mathrm{ip}^{([2]^{},X)}
=\eta^{(0,X)};$
\item[(iii)] $\mathrm{ip}^{([2]^{},0)\sharp}\circ\eta^{(0,X)}=
\mathrm{ip}^{([2]^{},X)}.
$
\end{itemize}

The reader is advised to consult the diagram appearing in Figure~\ref{FDCHDZ}.
\end{proposition}

\begin{figure}
\begin{center}
\begin{tikzpicture}
[ACliment/.style={-{To [angle'=45, length=5.75pt, width=4pt, round]}},scale=1]
\node[] (x) at (0,0) [] {$X$};
\node[] (pt) at (6,0) [] {$\mathrm{T}_{\Sigma}(X)$};
\node[] (2pq) at (6,-3) []  {$[\mathrm{Pth}_{\boldsymbol{\mathcal{A}}^{(2)}}]
$};

\draw[ACliment]  (x) 	to node [above right]	
{$\eta^{(0,X)}$} (pt);
\draw[ACliment, bend right=10]  (x) 	to node [below left]	
{$\mathrm{ip}^{([2]^{},X)}$} (2pq);

\node[] (B1) at (6,-1.5)  [] {};
\draw[ACliment]  ($(B1)+(0,1.2)$) to node [above, fill=white] {
$\mathrm{ip}^{([2]^{},0)\sharp}$
} ($(B1)+(0,-1.2)$);
\draw[ACliment, bend right]  ($(B1)+(.5,-1.2)$) to node [ below, fill=white] {
$ \mathrm{tg}^{(0,[2]^{})}$
} ($(B1)+(.5,1.2)$);
\draw[ACliment, bend left]  ($(B1)+(-.5,-1.2)$) to node [below, fill=white] {
$ \mathrm{sc}^{(0,[2]^{})}$
} ($(B1)+(-.5,1.2)$);

\end{tikzpicture}
\end{center}
\caption{Quotient mappings relative to $X$ at layers 0 \& 2.}
\label{FDCHDZ}
\end{figure}

\section{
\texorpdfstring
{A structure of partial  $\Sigma^{\boldsymbol{\mathcal{A}}^{(2)}}$-algebra on $[\mathrm{Pth}_{\boldsymbol{\mathcal{A}}^{(2)}}]$}
{An algebra on the quotient of second-order paths}
}

In this section, as a consequence of the main result of the previous chapter, we prove that the quotient $[\mathrm{Pth}_{\boldsymbol{\mathcal{A}}^{(2)}}] = \mathrm{Pth}_{\boldsymbol{\mathcal{A}}^{(2)}}/\mathrm{Ker}(\mathrm{CH}^{(2)})$ is equipped with a structure of partial $\Sigma^{\boldsymbol{\mathcal{A}}^{(2)}}$-algebra, we denote by $[\mathbf{Pth}_{\boldsymbol{\mathcal{A}}^{(2)}}]$ the corresponding $\Sigma^{\boldsymbol{\mathcal{A}}^{(2)}}$-algebra. Moreover, we consider several reducts of $[\mathbf{Pth}_{\boldsymbol{\mathcal{A}}^{(2)}}]$ and investigate some morphisms from and to these reducts of $[\mathbf{Pth}_{\boldsymbol{\mathcal{A}}^{(2)}}]$.

\begin{restatable}{proposition}{PDCHDCatAlg}
\label{PDCHDCatAlg}
\index{path!second-order!$[\mathbf{Pth}_{\boldsymbol{\mathcal{A}}^{(2)}}]$}
The $S$-sorted set $[\mathrm{Pth}_{\boldsymbol{\mathcal{A}}^{(2)}}]$ is equipped, in a natural way, with a structure of partial $\Sigma^{\boldsymbol{\mathcal{A}}^{(2)}}$-algebra. We denote by $[\mathbf{Pth}_{\boldsymbol{\mathcal{A}}^{(2)}}]$ the corresponding $\Sigma^{\boldsymbol{\mathcal{A}}^{(2)}}$-algebra (which is a quotient of $\mathbf{Pth}_{\boldsymbol{\mathcal{A}}^{(2)}}$, the partial $\Sigma^{\boldsymbol{\mathcal{A}}^{(2)}}$-algebra constructed in Proposition~\ref{PDPthDCatAlg}).
Moreover, the mapping $\mathrm{pr}^{\mathrm{Ker}(\mathrm{CH}^{(2)})}$, i.e.,
$$
\mathrm{pr}^{\mathrm{Ker}(\mathrm{CH}^{(2)})}
\colon
\mathbf{Pth}_{\boldsymbol{\mathcal{A}}^{(2)}}
\mor
[\mathbf{Pth}_{\boldsymbol{\mathcal{A}}^{(2)}}]
$$
is a closed and surjective $\Sigma^{\boldsymbol{\mathcal{A}}^{(2)}}$-homomorphism from $\mathbf{Pth}_{\boldsymbol{\mathcal{A}}^{(2)}}$ to $[\mathbf{Pth}_{\boldsymbol{\mathcal{A}}^{(2)}}]$.
\end{restatable}

\begin{remark}\label{RDCHMono} The monomorphic second-order Curry-Howard mapping, that is, 
$$
\mathrm{CH}^{(2)\mathrm{m}}
\colon
[\mathrm{Pth}_{\boldsymbol{\mathcal{A}}^{(2)}}]
\mor
\mathrm{T}_{\Sigma^{\boldsymbol{\mathcal{A}}^{(2)}}}(X)
$$
is not a $\Sigma^{\boldsymbol{\mathcal{A}}^{(2)}}$-homomorphism from 
$[\mathbf{Pth}_{\boldsymbol{\mathcal{A}}^{(2)}}]$ to 
$\mathbf{T}_{\Sigma^{\boldsymbol{\mathcal{A}}^{(2)}}}(X)$. The reader is advised to consult the counterexamples stated in Propositions~\ref{PDCHNotHomCat} and~\ref{PDCHNotHomDCat}. 
This will be rectified later by considering a suitable quotient of the free $\Sigma^{\boldsymbol{\mathcal{A}}^{(2)}}$-algebra on $X$.
\end{remark}

We next consider the $\Sigma^{\boldsymbol{\mathcal{A}}}$-reduct of the many-sorted partial $\Sigma^{\boldsymbol{\mathcal{A}}^{(2)}}$-algebra of second-order path classes and study its connections with the mappings presented in the previous section.

\begin{definition} For the partial $\Sigma^{\boldsymbol{\mathcal{A}}^{(2)}}$-algebra $[\mathbf{Pth}_{\boldsymbol{\mathcal{A}}^{(2)}}]$, we denote by $[\mathbf{Pth}^{(1,2)}_{\boldsymbol{\mathcal{A}}^{(2)}}]$ the partial $\Sigma^{\boldsymbol{\mathcal{A}}}$-algebra $\mathbf{in}^{\Sigma,(1,2)}\left([
\mathbf{Pth}_{\boldsymbol{\mathcal{A}}^{(2)}}]\right)$. We will call  $[\mathbf{Pth}^{(1,2)}_{\boldsymbol{\mathcal{A}}^{(2)}}]$ the $\Sigma^{\boldsymbol{\mathcal{A}}}$-reduct of the partial $\Sigma^{\boldsymbol{\mathcal{A}}^{(2)}}$-algebra  $[\mathbf{Pth}_{\boldsymbol{\mathcal{A}}^{(2)}}]$.
\end{definition}

\begin{proposition} The mapping $\mathrm{pr}^{\mathrm{Ker}(\mathrm{CH}^{(2)})}$
is a closed and surjective $\Sigma^{\boldsymbol{\mathcal{A}}}$-homomorphism from $\mathbf{Pth}^{(1,2)}_{\boldsymbol{\mathcal{A}}^{(2)}}$ to $[\mathbf{Pth}^{(1,2)}_{\boldsymbol{\mathcal{A}}}]$.
\end{proposition}

\begin{restatable}{proposition}{PDCHDU}
\label{PDCHDU} The many-sorted mappings $\mathrm{sc}^{([1],[2]^{})}$ and $\mathrm{tg}^{([1],[2]^{})}$ are closed $\Sigma^{\boldsymbol{\mathcal{A}}}$-homomorphisms from $[\mathbf{Pth}^{(1,2)}_{\boldsymbol{\mathcal{A}}^{(2)}}]$ to $[\mathbf{PT}_{\boldsymbol{\mathcal{A}}}]$.
\end{restatable}

\begin{restatable}{proposition}{PDCHDUIp}
\label{PDCHDUIp} The many-sorted mapping $\mathrm{ip}^{([2]^{},[1])\sharp}$  is a closed $\Sigma^{\boldsymbol{\mathcal{A}}}$-homomorphism from  $[\mathbf{PT}_{\boldsymbol{\mathcal{A}}}]$ to $[\mathbf{Pth}^{(1,2)}_{\boldsymbol{\mathcal{A}}^{(2)}}]$.
\end{restatable}

We next consider the $\Sigma$-reduct of the many-sorted partial $\Sigma^{\boldsymbol{\mathcal{A}}^{(2)}}$-algebra of second-order path first-order classes and study its connections with the mappings presented in the previous section.

\begin{definition} For the partial $\Sigma^{\boldsymbol{\mathcal{A}}^{(2)}}$-algebra $[\mathbf{Pth}_{\boldsymbol{\mathcal{A}}^{(2)}}]$, we denote by $[\mathbf{Pth}^{(0,2)}_{\boldsymbol{\mathcal{A}}^{(2)}}]$ the $\Sigma$-algebra $\mathbf{in}^{\Sigma,(0,2)}\left(
[\mathbf{Pth}_{\boldsymbol{\mathcal{A}}^{(2)}}]\right)$. We will call  $[\mathbf{Pth}^{(0,2)}_{\boldsymbol{\mathcal{A}}^{(2)}}]$ the $\Sigma$-reduct of the partial $\Sigma^{\boldsymbol{\mathcal{A}}^{(2)}}$-algebra  $[\mathbf{Pth}_{\boldsymbol{\mathcal{A}}^{(2)}}]$.
\end{definition}

\begin{proposition} The mapping $\mathrm{pr}^{\mathrm{Ker}(\mathrm{CH}^{(2)})}$
is a surjective $\Sigma$-homomorphism from $\mathbf{Pth}^{(0,2)}_{\boldsymbol{\mathcal{A}}^{(2)}}$ to $[\mathbf{Pth}^{(0,2)}_{\boldsymbol{\mathcal{A}}}]$.
\end{proposition}

\begin{restatable}{proposition}{PDCHDZ}
\label{PDCHDZ} The many-sorted mappings $\mathrm{sc}^{(0,[2]^{})}$ and $\mathrm{tg}^{(0,[2]^{})}$ are $\Sigma$-homomorphisms from $[\mathbf{Pth}^{(0,2)}_{\boldsymbol{\mathcal{A}}^{(2)}}]$ to $\mathbf{T}_{\Sigma}(X)$.
\end{restatable}

\begin{restatable}{proposition}{PDCHDZIp}
\label{PDCHDZIp} The many-sorted mapping $\mathrm{ip}^{([2]^{},0)\sharp}$  is a $\Sigma$-homomorphism from  $\mathbf{T}_{\Sigma}(X)$ to $[\mathbf{Pth}^{(0,2)}_{\boldsymbol{\mathcal{A}}^{(2)}}]$.
\end{restatable}

\section{
\texorpdfstring
{An Artinian order on $\coprod[\mathrm{Pth}_{\boldsymbol{\mathcal{A}}^{(2)}}]$}
{An Artinian order on the quotient of second-order paths}
}

In this section we define an order on $\coprod[\mathrm{Pth}_{\boldsymbol{\mathcal{A}}^{(2)}}]$, the coproduct of the quotient $[\mathrm{Pth}_{\boldsymbol{\mathcal{A}}^{(2)}}]$, formed by all labelled second-order path classes $([\mathfrak{P}^{(2)}]^{}_{s},s)$ with $s\in S$ and $[\mathfrak{P}^{(2)}]^{}_{s}\in[\mathrm{Pth}_{\boldsymbol{\mathcal{A}}^{(2)}}]_{}$.

\begin{restatable}{definition}{DDCHOrd}
\label{DDCHOrd} 
\index{partial order!second-order!$\leq_{[\mathbf{Pth}_{\boldsymbol{\mathcal{A}}^{(2)}}]}$}
Let $\leq_{[\mathbf{Pth}_{\boldsymbol{\mathcal{A}}^{(2)}}]}$ be the binary relation defined on $\coprod[\mathrm{Pth}_{\boldsymbol{\mathcal{A}}^{(2)}}]$ containing every pair
$(([\mathfrak{Q}^{(2)}]^{}_{t},t),([\mathfrak{P}^{(2)}]^{}_{s},s))$ in $(\coprod[\mathrm{Pth}_{\boldsymbol{\mathcal{A}}^{(2)}}])^{2}$ satisfying that
$$
\exists\, \mathfrak{Q}'^{(2)}\in 
\left[\mathfrak{Q}^{(2)}\right]^{}_{t},\,\, 
\exists\, \mathfrak{P}'^{(2)}\in \left[
\mathfrak{P}^{(2)}\right]^{}_{s}\, 
\left(
\left(
\mathfrak{Q}'^{(2)},t
\right)
\leq_{\mathbf{Pth}_{\boldsymbol{\mathcal{A}}^{(2)}}}
\left(\mathfrak{P}'^{(2)},s
\right)
\right).
$$
That is, $([\mathfrak{Q}^{(2)}]^{}_{t},t)$ $\leq_{[\mathbf{Pth}_{\boldsymbol{\mathcal{A}}^{(2)}}]}$-precedes $([\mathfrak{P}^{(2)}]^{}_{s},s)$ if there exists a pair of representative second-order paths $\mathfrak{Q}'^{(2)}\in [\mathfrak{Q}^{(2)}]^{}_{t}$ and $\mathfrak{P}'^{(2)}\in [\mathfrak{P}^{(2)}]^{}_{s}$ for which $(\mathfrak{Q}'^{(2)},t)$ $\leq_{\mathbf{Pth}_{\boldsymbol{\mathcal{A}}^{(2)}}}$-precedes $(\mathfrak{P}'^{(2)},s)$.
\end{restatable}

Our aim is to prove that $\leq_{[\mathbf{Pth}_{\boldsymbol{\mathcal{A}}^{(2)}}]}$ is an Artinian order on $[\mathrm{Pth}_{\boldsymbol{\mathcal{A}}^{(2)}}]$ and for this we need to state the following lemmas, which are ultimately based on the presentation given in~\cite{GK15}.

In the following lemma we show that, if we are given two second-order paths in the same class under the Kernel of the second-order Curry Howard mapping and one of these second-order paths is smaller than another second-order path with respect to the partial order $\leq_{\mathbf{Pth}_{\boldsymbol{\mathcal{A}}^{(2)}}}$, then we can complete this inequality on the other side with an equivalent second-order path.

\begin{lemma}\label{LDCHOrdI} Let $t,s$ be sorts in $S$, let $\mathfrak{Q}'^{(2)},\mathfrak{Q}^{(2)}$ be second-order paths in $\mathrm{Pth}_{\boldsymbol{\mathcal{A}}^{(2)},t}$ and let $\mathfrak{P}^{(2)}$ be a second-order path in $\mathrm{Pth}_{\boldsymbol{\mathcal{A}}^{(2)},s}$. If  $(\mathfrak{Q}^{(2)},t)\leq_{\mathbf{Pth}_{\boldsymbol{\mathcal{A}}^{(2)}}}(\mathfrak{P}^{(2)},s)$ and $\mathfrak{Q}'^{(2)}\in[\mathfrak{Q}^{(2)}]^{}_{t}$ then there exists a second-order path $\mathfrak{P}'^{(2)}\in [\mathfrak{P}^{(2)}]^{}_{s}$ satisfying that  
$$ \left(
\mathfrak{Q}'^{(2)},t
\right)\leq_{\mathbf{Pth}_{\boldsymbol{\mathcal{A}}^{(2)}}}
\left(
\mathfrak{P}'^{(2)},s
\right).$$
\end{lemma}

\begin{proof}
Let us recall from Remark~\ref{RDOrd} that $(\mathfrak{Q}^{(2)},t)\leq_{\mathbf{Pth}_{\boldsymbol{\mathcal{A}}^{(2)}}}(\mathfrak{P}^{(2)},s)$ if and only if $s=t$ and $\mathfrak{Q}^{(2)}=\mathfrak{P}^{(2)}$ or there exists a natural number $m\in\mathbb{N}-\{0\}$, a word $\mathbf{w}\in S^{\star}$ of length $\bb{\mathbf{w}}=m+1$, and a family of second-order paths $(\mathfrak{R}^{(2)}_{k})$ in $\mathrm{Pth}_{\boldsymbol{\mathcal{A}}^{(2)},\mathbf{w}}$ such that $w_{0}=t$, $\mathfrak{R}^{(2)}_{0}=\mathfrak{Q}^{(2)}$, $w_{m}=s$, $\mathfrak{R}^{(2)}_{m}=\mathfrak{P}^{(2)}$ and for every $k\in m$, $(\mathfrak{R}^{(2)}_{k}, w_{k})\prec_{\mathbf{Pth}_{\boldsymbol{\mathcal{A}}^{(2)}}} (\mathfrak{R}^{(2)}_{k+1}, w_{k+1})$.

The lemma holds trivially in case $s=t$ and $\mathfrak{Q}^{(2)}=\mathfrak{P}^{(2)}$. Therefore, it remains to prove the other case. 

We will prove the lemma by induction on $m\in\mathbb{N}-\{0\}$.

\textsf{Base step of the induction.}

For $m=1$ we have that $(\mathfrak{Q}^{(2)},t)\prec_{\mathbf{Pth}_{\boldsymbol{\mathcal{A}}^{(2)}}}(\mathfrak{P}^{(2)},s)$, hence following Definition~\ref{DDOrd} we are in one of the following situations
\begin{enumerate}
\item $\mathfrak{P}^{(2)}$ and $\mathfrak{Q}^{(2)}$ are $(2,[1])$-identity second-order paths of the form 
\begin{align*}
\mathfrak{P}^{(2)}&=
\mathrm{ip}^{(2,[1])\sharp}_{s}\left(
\left[
P
\right]_{s}
\right),
&
\mathfrak{Q}^{(2)}&=
\mathrm{ip}^{(2,[1])\sharp}_{t}\left(
\left[
Q\right]_{t}
\right),
\end{align*}
for some path term classes $[P]_{s}\in [\mathrm{PT}_{\boldsymbol{\mathcal{A}}}]_{s}$ and  $[Q]_{t}\in [\mathrm{PT}_{\boldsymbol{\mathcal{A}}}]_{t}$ and the following inequality holds
$$
\left(\left[Q\right]_{t},t\right)
<_{[\mathbf{PT}_{\boldsymbol{\mathcal{A}}}]
}
\left(\left[P\right]_{s},s\right),
$$
where $\leq_{[\mathbf{PT}_{\boldsymbol{\mathcal{A}}}]
}$ is the Artinian partial order on $\coprod [\mathrm{PT}_{\boldsymbol{\mathcal{A}}}]$ introduced in Definition~\ref{DPTQOrd}.
\item $\mathfrak{P}^{(2)}$ is a second-order path of length strictly greater than one containing at least one   second-order echelon, and if its first   second-order echelon occurs at position $i\in\bb{\mathfrak{P}^{(2)}}$, then
\subitem if $i=0$, then $\mathfrak{Q}^{(2)}$ is equal to $\mathfrak{P}^{(2),0,0}$ or $\mathfrak{P}^{(2),1,\bb{\mathfrak{P}^{(2)}}-1}$,
\subitem if $i>0$, then $\mathfrak{Q}^{(2)}$ is equal to $\mathfrak{P}^{(2),0,i-1}$ or $\mathfrak{P}^{(2),i,\bb{\mathfrak{P}^{(2)}}-1}$,
\item $\mathfrak{P}^{(2)}$ is an echelonless second-order path, then
\subitem if $\mathfrak{P}^{(2)}$ is not head-constant, then let $i\in\bb{\mathfrak{P}^{(2)}}$ be the  maximum index for which $\mathfrak{P}^{(2),0,i}$ is a head-constant second-order path, then $\mathfrak{Q}^{(2)}$ is equal to $\mathfrak{P}^{(2),0,i}$ or $\mathfrak{P}^{(2),i+1,\bb{\mathfrak{P}^{(2)}}-1}$,
\subitem if $\mathfrak{P}^{(2)}$ is head-constant but not coherent, then let $i\in\bb{\mathfrak{P}^{(2)}}$ be the  maximum index for which $\mathfrak{P}^{(2),0,i}$ is a coherent second-order path, then $\mathfrak{Q}^{(2)}$ is equal to $\mathfrak{P}^{(2),0,i}$ or $\mathfrak{P}^{(2),i+1,\bb{\mathfrak{P}^{(2)}}-1}$,
\subitem if $\mathfrak{P}^{(2)}$ is head-constant and coherent then $\mathfrak{Q}^{(2)}$ is one of the second-order paths we can extract from $\mathfrak{P}^{(2)}$ in virtue of Lemma~\ref{LDPthExtract}.
\end{enumerate}

If~(1), we note that both $\mathfrak{P}^{(2)}$ and $\mathfrak{Q}^{(2)}$ are $(2,[1])$-identity second-order paths. Then, if $\mathfrak{Q}'^{(2)}$ is a second-order path in $[\mathfrak{Q}^{(2)}]^{}_{t}$ then we conclude, in virtue  of Corollary~\ref{CDCHUId}, that $\mathfrak{Q}'^{(2)}=\mathfrak{Q}^{(2)}$. This case follows easily, since $\mathfrak{P}^{(2)}$ is a second-order path in $[\mathfrak{P}^{(2)}]^{}_{s}$ satisfying that  
$$ \left(
\mathfrak{Q}^{(2)},t
\right)\leq_{\mathbf{Pth}_{\boldsymbol{\mathcal{A}}^{(2)}}}
\left(
\mathfrak{P}^{(2)},s
\right).$$

This completes the Case~(1).

If~(2), that is, if $\mathfrak{P}^{(2)}$ is a second-order path of length strictly greater than one containing at least one   second-order echelon and its first   second-order echelon occurs at position $i\in\bb{\mathfrak{P}^{(2)}}$, we consider the case~(2.1) $i=0$ and the subcases~(2.1.1) $\mathfrak{Q}^{(2)}=\mathfrak{P}^{(2),0,0}$ and~(2.1.2) $\mathfrak{Q}^{(2)}=\mathfrak{P}^{(2),1,\bb{\mathfrak{P}^{(2)}}-1}$ and the case~(2.2) $i>0$ and the subcases~(2.2.1) $\mathfrak{Q}^{(2)}=\mathfrak{P}^{(2),0,i-1}$ and (2.2.2) $\mathfrak{Q}^{(2)}=\mathfrak{P}^{(2),i,\bb{\mathfrak{P}^{(2)}}-1}$.
Let us note that in each case, we have that $t=s$.

If~(2.1.1), i.e.,  if $\mathfrak{P}^{(2)}$ is a second-order path of length strictly greater than one containing its first   second-order echelon on its first step and $\mathfrak{Q}^{(2)}=\mathfrak{P}^{(2),0,0}$ then, following Definition~\ref{DDCH}, the value of the second-order Curry-Howard mapping at $\mathfrak{P}^{(2)}$ is given by
$$
\mathrm{CH}^{(2)}_{s}\left(
\mathfrak{P}^{(2)}\right)
=
\mathrm{CH}^{(2)}_{s}\left(
\mathfrak{P}^{(2),1,\bb{\mathfrak{P}^{(2)}}-1}
\right)
\circ^{1\mathbf{T}_{\Sigma^{\boldsymbol{\mathcal{A}}^{(2)}}}(X)}_{s}
\mathrm{CH}^{(2)}_{s}\left(
\mathfrak{Q}^{(2)}
\right).
$$

Since $\mathfrak{Q}^{(2)}$ is a   second-order echelon, then in virtue of Proposition~\ref{PDCHEch} we have that $\mathfrak{Q}'^{(2)}=\mathfrak{Q}^{(2)}$. This case follows easily, since $\mathfrak{P}^{(2)}$ is a second-order path in $[\mathfrak{P}^{(2)}]^{}_{s}$ satisfying that  
$$ \left(
\mathfrak{Q}^{(2)},t
\right)\leq_{\mathbf{Pth}_{\boldsymbol{\mathcal{A}}^{(2)}}}
\left(
\mathfrak{P}^{(2)},s
\right).$$

The case $i=0$ and $\mathfrak{Q}^{(2)}=\mathfrak{P}^{(2),0,0}$ follows.

If~(2.1.2), i.e.,  if $\mathfrak{P}^{(2)}$ is a second-order path of length strictly greater than one containing its first   second-order echelon on its first step and  $\mathfrak{Q}^{(2)}=\mathfrak{P}^{(2),1,\bb{\mathfrak{P}^{(2)}}-1}$ then, following Definition~\ref{DDCH}, the value of the second-order Curry-Howard mapping at $\mathfrak{P}^{(2)}$ is given by
$$
\mathrm{CH}^{(2)}_{s}\left(
\mathfrak{P}^{(2)}
\right)
=
\mathrm{CH}^{(2)}_{s}\left(
\mathfrak{Q}^{(2)}
\right)
\circ^{1\mathbf{T}_{\Sigma^{\boldsymbol{\mathcal{A}}^{(2)}}}(X)}_{s}
\mathrm{CH}^{(2)}_{s}\left(
\mathfrak{P}^{(2),0,0}
\right).
$$

 Thus, if $\mathfrak{Q}'^{(2)}$ is any-second-order path in $[\mathfrak{Q}^{(2)}]^{}_{s}$ then, in virtue of Lemma~\ref{LDCH}, we have that 
$
\mathrm{sc}^{([1],2)}_{s}(\mathfrak{Q}'^{(2)})
=
\mathrm{sc}^{([1],2)}_{s}(\mathfrak{Q}^{(2)})
$. Moreover, the length of $\mathfrak{Q}'^{(2)}$ equals the length of $\mathfrak{Q}^{(2)}$. Therefore, since $\mathfrak{P}^{(2)}$ decomposes as the $1$-composition $\mathfrak{P}^{(2)}=\mathfrak{Q}^{(2)}\circ^{1\mathbf{Pth}_{\boldsymbol{\mathcal{A}}^{(2)}}}_{s}\mathfrak{P}^{(2),0,0}$, we have that 
$
\mathrm{sc}^{([1],2)}_{s}(\mathfrak{Q}^{(2)})
=
\mathrm{tg}^{([1],2)}_{s}(\mathfrak{P}^{(2),0,0}).
$ All in all, we conclude that 
$
\mathrm{sc}^{([1],2)}_{s}(\mathfrak{Q}'^{(2)})
=
\mathrm{tg}^{([1],2)}_{s}(\mathfrak{P}^{(2),0,0}).
$ 
Therefore the second-order path $\mathfrak{P}'^{(2)}$ defined by 
$$
\mathfrak{P}'^{(2)}
=\mathfrak{Q}'^{(2)}\circ^{1\mathbf{Pth}_{\boldsymbol{\mathcal{A}}^{(2)}}}_{s}\mathfrak{P}^{(2),0,0}
$$
 is a second-order path of length strictly greater than one containing a   second-order echelon on its first step. Taking into account Definition~\ref{DDOrd}, we have that 
$
(\mathfrak{Q}'^{(2)},s)
\prec_{\mathbf{Pth}_{\boldsymbol{\mathcal{A}}^{(2)}}}
(\mathfrak{P}'^{(2)},s).
$ Moreover the value of the second-order Curry-Howard mapping at $\mathfrak{P}'^{(2)}$ is given by
\begin{align*}
\mathrm{CH}^{(2)}_{s}\left(
\mathfrak{P}'^{(2)}
\right)
&=
\mathrm{CH}^{(2)}_{s}\left(
\mathfrak{Q}'^{(2)}
\right)
\circ^{1\mathbf{T}_{\Sigma^{\boldsymbol{\mathcal{A}}^{(2)}}}(X)}_{s}
\mathrm{CH}^{(2)}_{s}\left(
\mathfrak{P}^{(2),0,0}
\right)
\\&=
\mathrm{CH}^{(2)}_{s}\left(
\mathfrak{Q}^{(2)}
\right)
\circ^{1\mathbf{T}_{\Sigma^{\boldsymbol{\mathcal{A}}^{(2)}}}(X)}_{s}
\mathrm{CH}^{(2)}_{s}\left(
\mathfrak{P}^{(2),0,0}
\right)
\\&=
\mathrm{CH}^{(2)}_{s}\left(
\mathfrak{P}^{(2)}
\right).
\end{align*}

The case $i=0$ and $\mathfrak{Q}^{(2)}=\mathfrak{P}^{(2),1,\bb{\mathfrak{P}^{(2)}}-1}$ follows.

If~(2.2.1), i.e.,  if $\mathfrak{P}^{(2)}$ is a second-order path of length strictly greater than one containing its first   second-order echelon on a step $i\in\bb{\mathfrak{P}^{(2)}}$ different from the initial one and 
 $\mathfrak{Q}^{(2)}=\mathfrak{P}^{(2),0,i-1}$ then, following Definition~\ref{DDCH}, the value of the second-order Curry-Howard mapping at $\mathfrak{P}^{(2)}$ is given by
$$
\mathrm{CH}^{(2)}_{s}\left(
\mathfrak{P}^{(2)}
\right)
=
\mathrm{CH}^{(2)}_{s}\left(
\mathfrak{P}^{(2),i,\bb{\mathfrak{P}^{(2)}}-1}
\right)
\circ^{1\mathbf{T}_{\Sigma^{\boldsymbol{\mathcal{A}}^{(2)}}}(X)}_{s}
\mathrm{CH}^{(2)}_{s}\left(
\mathfrak{Q}^{(2)}
\right).
$$

Since $\mathfrak{Q}^{(2)}$ is an echelonless second-order path, then in virtue of Lemma~\ref{LDCHNEch} we have that $$\mathrm{CH}^{(2)}_{s}\left(\mathfrak{Q}^{(2)}
\right)\in 
\mathrm{T}_{\Sigma^{\boldsymbol{\mathcal{A}}^{(2)}}}(X)_{s}
-\left(
\eta^{(2,1)\sharp}\left[
\mathrm{PT}_{\boldsymbol{\mathcal{A}}}
\right]_{s}
\cup
\eta^{(2,\mathcal{A}^{(2)})}\left[
\mathcal{A}^{(2)}
\right]^{\mathrm{pct}}_{s}
\right)
.$$
Thus, if $\mathfrak{Q}'^{(2)}$ is any-second-order path in $[\mathfrak{Q}^{(2)}]^{}_{s}$ then, in virtue of Lemma~\ref{LDCHNEch} we conclude that $\mathfrak{Q}'^{(2)}$ is an echelonless second-order path. Moreover, in virtue of Lemma~\ref{LDCH}, we have that 
$
\mathrm{tg}^{([1],2)}_{s}(\mathfrak{Q}'^{(2)})
=
\mathrm{tg}^{([1],2)}_{s}(\mathfrak{Q}^{(2)})
$ and the length of $\mathfrak{Q}'^{(2)}$ equals the length of $\mathfrak{Q}^{(2)}$. Therefore, since $\mathfrak{P}^{(2)}$ decomposes as the $1$-composition $\mathfrak{P}^{(2)}=\mathfrak{P}^{(2),i,\bb{\mathfrak{P}^{(2)}}-1}\circ^{1\mathbf{Pth}_{\boldsymbol{\mathcal{A}}^{(2)}}}_{s}\mathfrak{Q}^{(2)}$, we have that 
$
\mathrm{sc}^{([1],2)}_{s}(\mathfrak{P}^{(2),i,\bb{\mathfrak{P}^{(2)}}-1})
=
\mathrm{tg}^{([1],2)}_{s}(\mathfrak{Q}^{(2)}).
$ All in all, we conclude that 
$
\mathrm{sc}^{([1],2)}_{s}(\mathfrak{P}^{(2),i,\bb{\mathfrak{P}^{(2)}}-1})
=
\mathrm{tg}^{([1],2)}_{s}(\mathfrak{Q}'^{(2)}).
$ 
Therefore the second-order path $\mathfrak{P}'^{(2)}$ defined by 
$$
\mathfrak{P}'^{(2)}
=\mathfrak{P}^{(2),i,\bb{\mathfrak{P}^{(2)}}-1}\circ^{1\mathbf{Pth}_{\boldsymbol{\mathcal{A}}^{(2)}}}_{s}\mathfrak{Q}'^{(2)}
$$
 is a second-order path of length strictly greater than one containing its first   second-order echelon at position $i$. Taking into account Definition~\ref{DDOrd}, we have that
$
(\mathfrak{Q}'^{(2)},s)
\prec_{\mathbf{Pth}_{\boldsymbol{\mathcal{A}}^{(2)}}}
(\mathfrak{P}'^{(2)},s).
$ Moreover  the value of the second-order Curry-Howard mapping at $\mathfrak{P}'^{(2)}$ is given by
\begin{align*}
\mathrm{CH}^{(2)}_{s}\left(
\mathfrak{P}'^{(2)}
\right)
&=
\mathrm{CH}^{(2)}_{s}\left(
\mathfrak{P}^{(2),i,\bb{\mathfrak{P}^{(2)}}-1}
\right)
\circ^{1\mathbf{T}_{\Sigma^{\boldsymbol{\mathcal{A}}^{(2)}}}(X)}_{s}
\mathrm{CH}^{(2)}_{s}\left(
\mathfrak{Q}'^{(2)}
\right)
\\&=
\mathrm{CH}^{(2)}_{s}\left(
\mathfrak{P}^{(2),i,\bb{\mathfrak{P}^{(2)}}-1}
\right)
\circ^{1\mathbf{T}_{\Sigma^{\boldsymbol{\mathcal{A}}^{(2)}}}(X)}_{s}
\mathrm{CH}^{(2)}_{s}\left(
\mathfrak{Q}^{(2)}
\right)
\\&=
\mathrm{CH}^{(2)}_{s}\left(
\mathfrak{P}^{(2)}
\right).
\end{align*}

The case $i>0$ and $\mathfrak{Q}^{(2)}=\mathfrak{P}^{(2),0,i-1}$ follows.

If~(2.2.2), i.e.,  if $\mathfrak{P}^{(2)}$ is a second-order path of length strictly greater than one containing its first   second-order echelon on a step $i\in\bb{\mathfrak{P}^{(2)}}$ different from the initial one and 
 $\mathfrak{Q}^{(2)}=\mathfrak{P}^{(2),i,\bb{\mathfrak{P}^{(2)}}-1}$ then, following Definition~\ref{DDCH}, the value of the second-order Curry-Howard mapping at $\mathfrak{P}^{(2)}$ is given by
$$
\mathrm{CH}^{(2)}_{s}\left(
\mathfrak{P}^{(2)}
\right)
=
\mathrm{CH}^{(2)}_{s}\left(
\mathfrak{Q}^{(2)}
\right)
\circ^{1\mathbf{T}_{\Sigma^{\boldsymbol{\mathcal{A}}^{(2)}}}(X)}_{s}
\mathrm{CH}^{(2)}_{s}\left(
\mathfrak{P}^{(2),0,i-1}
\right).
$$

Since $\mathfrak{Q}^{(2)}$ is a second-order path of length at least one containing a   second-order echelon on its first step, then in virtue of Lemmas~\ref{LDCHDEch} and~\ref{LDCHEchInt} we have that $$\mathrm{CH}^{(2)}_{s}\left(\mathfrak{Q}^{(2)}
\right)\in 
\eta^{(2,\mathcal{A}^{(2)})}\left[\mathcal{A}^{(2)}
\right]_{s}
\cup
\eta^{(2,\mathcal{A}^{(2)})}\left[
\mathcal{A}^{(2)}
\right]^{\mathrm{int}}_{s}
.$$
Thus, if $\mathfrak{Q}'^{(2)}$ is any-second-order path in $[\mathfrak{Q}^{(2)}]^{}_{s}$ then, in virtue of Lemmas~\ref{LDCHDEch} or~\ref{LDCHEchInt} we conclude that $\mathfrak{Q}'^{(2)}$  is a second-order path of length at least one containing a   second-order echelon on its first position. Moreover, in virtue of Lemma~\ref{LDCH}, we have that 
$
\mathrm{sc}^{([1],2)}_{s}(\mathfrak{Q}'^{(2)})
=
\mathrm{sc}^{([1],2)}_{s}(\mathfrak{Q}^{(2)})
$ and the length of $\mathfrak{Q}'^{(2)}$ equals the length of $\mathfrak{Q}^{(2)}$. Therefore, since $\mathfrak{P}^{(2)}$ decomposes as the $1$-composition $\mathfrak{P}^{(2)}=\mathfrak{Q}^{(2)}\circ^{1\mathbf{Pth}_{\boldsymbol{\mathcal{A}}^{(2)}}}_{s}\mathfrak{P}^{(2),0,i-1}$, we have that 
$
\mathrm{sc}^{([1],2)}_{s}(\mathfrak{Q}^{(2)})
=
\mathrm{tg}^{([1],2)}_{s}(\mathfrak{P}^{(2),0,i-1}).
$ All in all, we conclude that 
$
\mathrm{sc}^{([1],2)}_{s}(\mathfrak{Q}'^{(2)})
=
\mathrm{tg}^{([1],2)}_{s}(\mathfrak{P}^{(2),0,i-1}).
$ 
Therefore the second-order path $\mathfrak{P}'^{(2)}$ defined by 
$$
\mathfrak{P}'^{(2)}
=\mathfrak{Q}'^{(2)}\circ^{1\mathbf{Pth}_{\boldsymbol{\mathcal{A}}^{(2)}}}_{s}\mathfrak{P}^{(2),0,i-1}
$$
 is a second-order path of length strictly greater than one containing its first   second-order echelon at position $i$. Taking into account Definition~\ref{DDOrd}, we have that
$
(\mathfrak{Q}'^{(2)},s)
\prec_{\mathbf{Pth}_{\boldsymbol{\mathcal{A}}^{(2)}}}
(\mathfrak{P}'^{(2)},s).
$ Moreover the value of the second-order Curry-Howard mapping at $\mathfrak{P}'^{(2)}$ is given by
\begin{align*}
\mathrm{CH}^{(2)}_{s}\left(
\mathfrak{P}'^{(2)}
\right)
&=
\mathrm{CH}^{(2)}_{s}\left(
\mathfrak{Q}'^{(2)}
\right)
\circ^{1\mathbf{T}_{\Sigma^{\boldsymbol{\mathcal{A}}^{(2)}}}(X)}_{s}
\mathrm{CH}^{(2)}_{s}\left(
\mathfrak{P}^{(2),0,i-1}
\right)
\\&=
\mathrm{CH}^{(2)}_{s}\left(
\mathfrak{Q}^{(2)}
\right)
\circ^{1\mathbf{T}_{\Sigma^{\boldsymbol{\mathcal{A}}^{(2)}}}(X)}_{s}
\mathrm{CH}^{(2)}_{s}\left(
\mathfrak{P}^{(2),0,i-1}
\right)
\\&=
\mathrm{CH}^{(2)}_{s}\left(\mathfrak{P}^{(2)}
\right).
\end{align*}

The case $i>0$ and $\mathfrak{Q}^{(2)}=\mathfrak{P}^{(2),i,\bb{\mathfrak{P}^{(2)}}-1}$ follows.

This completes the Case~(2).
 
If~(3), that is, if $\mathfrak{P}^{(2)}$ is an echelonless second-order path, then we consider different cases according to the nature of $\mathfrak{P}^{(2)}$. In this regard, we consider the case (3.1), where $\mathfrak{P}^{(2)}$ is an echelonless second-order path that is not head-constant, the case (3.2), where $\mathfrak{P}^{(2)}$ is a head-constant echelonless second-order path that is not coherent, and the case~(3.3), where $\mathfrak{P}^{(2)}$ is a coherent head-constant echelonless second-order path.

If~(3.1), i.e., if $\mathfrak{P}^{(2)}$ is an echelonless second-order path that is not head-constant, then let $i\in\bb{\mathfrak{P}^{(2)}}$ be the maximum index for which $\mathfrak{P}^{(2),0,i}$ is a head-constant echelonless second-order path. Then we have two possibilities for $\mathfrak{Q}^{(2)}$, subcase~(3.1.1) $\mathfrak{Q}^{(2)}=\mathfrak{P}^{(2),0,i}$ and subcase~(3.1.2) $\mathfrak{Q}^{(2)}=\mathfrak{P}^{(2),i+1,\bb{\mathfrak{P}^{(2)}}-1}$.

If~(3.1.1), i.e., if $\mathfrak{P}^{(2)}$ is an echelonless second-order path that is not head-constant, $i\in\bb{\mathfrak{P}^{(2)}}$ is the maximum index for which $\mathfrak{P}^{(2),0,i}$ is a head-constant echelonless second-order path and $\mathfrak{Q}^{(2)}=\mathfrak{P}^{(2),0,i}$ then, following Definition~\ref{DDCH}, the value of the second-order Curry-Howard mapping at $\mathfrak{P}^{(2)}$ is given by
$$
\mathrm{CH}^{(2)}_{s}\left(
\mathfrak{P}^{(2)}
\right)
=
\mathrm{CH}^{(2)}_{s}\left(
\mathfrak{P}^{(2),i+1,\bb{\mathfrak{P}^{(2)}}-1}
\right)
\circ^{1\mathbf{T}_{\Sigma^{\boldsymbol{\mathcal{A}}^{(2)}}}(X)}_{s}
\mathrm{CH}^{(2)}_{s}\left(
\mathfrak{Q}^{(2)}
\right).
$$

Since $\mathfrak{Q}^{(2)}$ is a head-constant echelonless second-order path, then in virtue of Lemma~\ref{LDCHNEchHd} we have that $\mathrm{CH}^{(2)}_{s}(\mathfrak{Q}^{(2)})\in [\mathrm{T}_{\Sigma^{\boldsymbol{\mathcal{A}}^{(2)}}}(X)]^{\mathsf{HdC}}_{s}
$. 
Thus, if $\mathfrak{Q}'^{(2)}$ is any-second-order path in $[\mathfrak{Q}^{(2)}]^{}_{s}$ then, in virtue of Lemma~\ref{LDCHNEchHd} we conclude that $\mathfrak{Q}'^{(2)}$ is a head-constant echelonless second-order path of length at least one. Moreover, in virtue of Lemma~\ref{LDCH}, we have that 
$
\mathrm{tg}^{([1],2)}_{s}(\mathfrak{Q}'^{(2)})
=
\mathrm{tg}^{([1],2)}_{s}(\mathfrak{Q}^{(2)})
$ and the length of $\mathfrak{Q}'^{(2)}$ equals the length of $\mathfrak{Q}^{(2)}$. Therefore, since $\mathfrak{P}^{(2)}$ decomposes as the $1$-composition $\mathfrak{P}^{(2)}=\mathfrak{P}^{(2),i+1,\bb{\mathfrak{P}^{(2)}}-1}\circ^{1\mathbf{Pth}_{\boldsymbol{\mathcal{A}}^{(2)}}}_{s}\mathfrak{Q}^{(2)}$, we have that 
$
\mathrm{sc}^{([1],2)}_{s}(\mathfrak{P}^{(2),i+1,\bb{\mathfrak{P}^{(2)}}-1})
=
\mathrm{tg}^{([1],2)}_{s}(\mathfrak{Q}^{(2)}).
$ All in all, we conclude that 
$
\mathrm{sc}^{([1],2)}_{s}(\mathfrak{P}^{(2),i+1,\bb{\mathfrak{P}^{(2)}}-1})
=
\mathrm{tg}^{([1],2)}_{s}(\mathfrak{Q}'^{(2)}).
$ 
Therefore the second-order path $\mathfrak{P}'^{(2)}$ defined by 
$$
\mathfrak{P}'^{(2)}
=\mathfrak{P}^{(2),i+1,\bb{\mathfrak{P}^{(2)}}-1}\circ^{1\mathbf{Pth}_{\boldsymbol{\mathcal{A}}^{(2)}}}_{s}\mathfrak{Q}'^{(2)}
$$
 is an echelonless second-order path that is not head-constant for which $i\in\bb{\mathfrak{P}^{(2)}}$ is the maximum index for which $\mathfrak{P}'^{(2),0,i}=\mathfrak{Q}'^{(2)}$ is a head-constant second-order path. Taking into account Definition~\ref{DDOrd}, we have that 
$
(\mathfrak{Q}'^{(2)},s)
\prec_{\mathbf{Pth}_{\boldsymbol{\mathcal{A}}^{(2)}}}
(\mathfrak{P}'^{(2)},s).
$ Moreover the value of the second-order Curry-Howard mapping at $\mathfrak{P}'^{(2)}$ is given by
\begin{align*}
\mathrm{CH}^{(2)}_{s}\left(
\mathfrak{P}'^{(2)}
\right)
&=
\mathrm{CH}^{(2)}_{s}\left(
\mathfrak{P}^{(2),i+1,\bb{\mathfrak{P}^{(2)}}-1}
\right)
\circ^{1\mathbf{T}_{\Sigma^{\boldsymbol{\mathcal{A}}^{(2)}}}(X)}_{s}
\mathrm{CH}^{(2)}_{s}\left(
\mathfrak{Q}'^{(2)}
\right)
\\&=
\mathrm{CH}^{(2)}_{s}\left(
\mathfrak{P}^{(2),i+1,\bb{\mathfrak{P}^{(2)}}-1}
\right)
\circ^{1\mathbf{T}_{\Sigma^{\boldsymbol{\mathcal{A}}^{(2)}}}(X)}_{s}
\mathrm{CH}^{(2)}_{s}\left(
\mathfrak{Q}^{(2)}
\right)
\\&=
\mathrm{CH}^{(2)}_{s}\left(\mathfrak{P}^{(2)}
\right).
\end{align*}

The case $\mathfrak{Q}^{(2)}=\mathfrak{P}^{(2),0,i}$ follows.

If~(3.1.2), i.e., if $\mathfrak{P}^{(2)}$ is an echelonless second-order path that is not head-constant, $i\in\bb{\mathfrak{P}^{(2)}}$ is the maximum index for which $\mathfrak{P}^{(2),0,i}$ is a head-constant echelonless second-order path and $\mathfrak{Q}^{(2)}=\mathfrak{P}^{(2),i+1,\bb{\mathfrak{P}^{(2)}}-1}$ then, following Definition~\ref{DDCH}, the value of the second-order Curry-Howard mapping at $\mathfrak{P}^{(2)}$ is given by
$$
\mathrm{CH}^{(2)}_{s}\left(
\mathfrak{P}^{(2)}
\right)
=
\mathrm{CH}^{(2)}_{s}\left(
\mathfrak{Q}^{(2)}
\right)
\circ^{1\mathbf{T}_{\Sigma^{\boldsymbol{\mathcal{A}}^{(2)}}}(X)}_{s}
\mathrm{CH}^{(2)}_{s}\left(
\mathfrak{P}^{(2),0,i}
\right).
$$

Since $\mathfrak{Q}^{(2)}$ is an echelonless second-order path, then in virtue of Lemma~\ref{LDCHNEch} we have that $$\mathrm{CH}^{(2)}_{s}\left(\mathfrak{Q}^{(2)}
\right)\in 
\mathrm{T}_{\Sigma^{\boldsymbol{\mathcal{A}}^{(2)}}}(X)_{s}
-
\left(
\eta^{(2,1)\sharp}
\left[\mathrm{PT}_{\boldsymbol{\mathcal{A}}}
\right]_{s}
\cup
\eta^{(2,\mathcal{A}^{(2)})}
\left[\mathcal{A}^{(2)}
\right]_{s}
\right)
.$$
Thus, if $\mathfrak{Q}'^{(2)}$ is any-second-order path in $[\mathfrak{Q}^{(2)}]^{}_{s}$ then, in virtue of Lemma~\ref{LDCHNEch} we conclude that $\mathfrak{Q}'^{(2)}$  is an echelonless second-order path. Moreover, in virtue of Lemma~\ref{LDCH}, we have that 
$
\mathrm{sc}^{([1],2)}_{s}(\mathfrak{Q}'^{(2)})
=
\mathrm{sc}^{([1],2)}_{s}(\mathfrak{Q}^{(2)})
$ and the length of $\mathfrak{Q}'^{(2)}$ equals the length of $\mathfrak{Q}^{(2)}$. Therefore, since $\mathfrak{P}^{(2)}$ decomposes as the $1$-composition $\mathfrak{P}^{(2)}=\mathfrak{Q}^{(2)}\circ^{1\mathbf{Pth}_{\boldsymbol{\mathcal{A}}^{(2)}}}_{s}\mathfrak{P}^{(2),0,i}$, we have that 
$
\mathrm{sc}^{([1],2)}_{s}(\mathfrak{Q}^{(2)})
=
\mathrm{tg}^{([1],2)}_{s}(\mathfrak{P}^{(2),0,i}).
$ All in all, we conclude that 
$
\mathrm{sc}^{([1],2)}_{s}(\mathfrak{Q}'^{(2)})
=
\mathrm{tg}^{([1],2)}_{s}(\mathfrak{P}^{(2),0,i}).
$ 
Therefore the second-order path $\mathfrak{P}'^{(2)}$ defined by 
$$
\mathfrak{P}'^{(2)}
=\mathfrak{Q}'^{(2)}\circ^{1\mathbf{Pth}_{\boldsymbol{\mathcal{A}}^{(2)}}}_{s}\mathfrak{P}^{(2),0,i}
$$
 is an echelonless second-order path that is not head-constant for which $i\in\bb{\mathfrak{P}^{(2)}}$ is the maximum index for which $\mathfrak{P}'^{(2),0,i}=\mathfrak{P}^{(2),0,i}$ is a head-constant second-order path.  Taking into account Definition~\ref{DDOrd}, we have that
$
(\mathfrak{Q}'^{(2)},s)
\prec_{\mathbf{Pth}_{\boldsymbol{\mathcal{A}}^{(2)}}}
(\mathfrak{P}'^{(2)},s).
$ Moreover the value of the second-order Curry-Howard mapping at $\mathfrak{P}'^{(2)}$ is given by
\begin{align*}
\mathrm{CH}^{(2)}_{s}\left(
\mathfrak{P}'^{(2)}
\right)
&=
\mathrm{CH}^{(2)}_{s}\left(
\mathfrak{Q}'^{(2)}
\right)
\circ^{1\mathbf{T}_{\Sigma^{\boldsymbol{\mathcal{A}}^{(2)}}}(X)}_{s}
\mathrm{CH}^{(2)}_{s}\left(
\mathfrak{P}^{(2),0,i}
\right)
\\&=
\mathrm{CH}^{(2)}_{s}\left(
\mathfrak{Q}^{(2)}
\right)
\circ^{1\mathbf{T}_{\Sigma^{\boldsymbol{\mathcal{A}}^{(2)}}}(X)}_{s}
\mathrm{CH}^{(2)}_{s}\left(
\mathfrak{P}^{(2),0,i}
\right)
\\&=
\mathrm{CH}^{(2)}_{s}
\left(
\mathfrak{P}^{(2)}
\right).
\end{align*}

The case $\mathfrak{Q}^{(2)}=\mathfrak{P}^{(2),i+1,\bb{\mathfrak{P}^{(2)}}-1}$ follows.

This completes the case (3.1) of $\mathfrak{P}^{(2)}$ being an echelonless second-order path that is not head-constant. 

We move to case (3.2), where $\mathfrak{P}^{(2)}$ is a head-constant echelonless second-order path that is not coherent. Let $i\in\bb{\mathfrak{P}^{(2)}}$ be the maximum index for which $\mathfrak{P}^{(2),0,i}$ is a coherent head-constant echelonless second-order path. Then we have two possibilities for $\mathfrak{Q}^{(2)}$, subcase~(3.2.1) $\mathfrak{Q}^{(2)}=\mathfrak{P}^{(2),0,i}$ and subcase~(3.2.2) $\mathfrak{Q}^{(2)}=\mathfrak{P}^{(2),i+1,\bb{\mathfrak{P}^{(2)}}-1}$.

If~(3.2.1), i.e., if $\mathfrak{P}^{(2)}$ is a head-constant echelonless second-order path that is not coherent, $i\in\bb{\mathfrak{P}^{(2)}}$ is the maximum index for which $\mathfrak{P}^{(2),0,i}$ is a coherent head-constant echelonless second-order path and $\mathfrak{Q}^{(2)}=\mathfrak{P}^{(2),0,i}$ then, following Definition~\ref{DDCH}, the value of the second-order Curry-Howard mapping at $\mathfrak{P}^{(2)}$ is given by
$$
\mathrm{CH}^{(2)}_{s}\left(
\mathfrak{P}^{(2)}
\right)
=
\mathrm{CH}^{(2)}_{s}\left(
\mathfrak{P}^{(2),i+1,\bb{\mathfrak{P}^{(2)}}-1}
\right)
\circ^{1\mathbf{T}_{\Sigma^{\boldsymbol{\mathcal{A}}^{(2)}}}(X)}_{s}
\mathrm{CH}^{(2)}_{s}\left(
\mathfrak{Q}^{(2)}
\right).
$$

Since $\mathfrak{Q}^{(2)}$ is a coherent head-constant echelonless second-order path, then in virtue of Lemma~\ref{LDCHNEchHdC} we have that $\mathrm{CH}^{(2)}_{s}(\mathfrak{Q}^{(2)})\in [\mathrm{T}_{\Sigma^{\boldsymbol{\mathcal{A}}^{(2)}}}(X)]^{\mathsf{HdC}\And\mathsf{C}}_{s}
$. 
Thus, if $\mathfrak{Q}'^{(2)}$ is any-second-order path in $[\mathfrak{Q}^{(2)}]^{}_{s}$ then, in virtue of Lemma~\ref{LDCHNEchHdC} we conclude that $\mathfrak{Q}'^{(2)}$ is a coherent head-constant echelonless second-order path. Moreover, in virtue of Lemma~\ref{LDCH}, we have that 
$
\mathrm{tg}^{([1],2)}_{s}(\mathfrak{Q}'^{(2)})
=
\mathrm{tg}^{([1],2)}_{s}(\mathfrak{Q}^{(2)})
$ and the length of $\mathfrak{Q}'^{(2)}$ equals the length of $\mathfrak{Q}^{(2)}$. Therefore, since $\mathfrak{P}^{(2)}$ decomposes as the $1$-composition $\mathfrak{P}^{(2)}=\mathfrak{P}^{(2),i+1,\bb{\mathfrak{P}^{(2)}}-1}\circ^{1\mathbf{Pth}_{\boldsymbol{\mathcal{A}}^{(2)}}}_{s}\mathfrak{Q}^{(2)}$, we have that 
$
\mathrm{sc}^{([1],2)}_{s}(\mathfrak{P}^{(2),i+1,\bb{\mathfrak{P}^{(2)}}-1})
=
\mathrm{tg}^{([1],2)}_{s}(\mathfrak{Q}^{(2)}).
$ All in all, we conclude that 
$
\mathrm{sc}^{([1],2)}_{s}(\mathfrak{P}^{(2),i+1,\bb{\mathfrak{P}^{(2)}}-1})
=
\mathrm{tg}^{([1],2)}_{s}(\mathfrak{Q}'^{(2)}).
$ 
Therefore the second-order path $\mathfrak{P}'^{(2)}$ defined by 
$$
\mathfrak{P}'^{(2)}
=\mathfrak{P}^{(2),i+1,\bb{\mathfrak{P}^{(2)}}-1}\circ^{1\mathbf{Pth}_{\boldsymbol{\mathcal{A}}^{(2)}}}_{s}\mathfrak{Q}'^{(2)}
$$
is, in virtue of Lemma~\ref{LTech}, a head-constant echelonless second-order path that is not coherent for which $i\in\bb{\mathfrak{P}^{(2)}}$ is the maximum index for which $\mathfrak{P}'^{(2),0,i}=\mathfrak{Q}'^{(2)}$ is a coherent head-constant echelonless second-order path. Taking into account Definition~\ref{DDOrd}, we have that 
$
(\mathfrak{Q}'^{(2)},s)
\prec_{\mathbf{Pth}_{\boldsymbol{\mathcal{A}}^{(2)}}}
(\mathfrak{P}'^{(2)},s).
$ Moreover the value of the second-order Curry-Howard mapping at $\mathfrak{P}'^{(2)}$ is given by
\begin{align*}
\mathrm{CH}^{(2)}_{s}\left(
\mathfrak{P}'^{(2)}
\right)
&=
\mathrm{CH}^{(2)}_{s}\left(
\mathfrak{P}^{(2),i+1,\bb{\mathfrak{P}^{(2)}}-1}
\right)
\circ^{1\mathbf{T}_{\Sigma^{\boldsymbol{\mathcal{A}}^{(2)}}}(X)}_{s}
\mathrm{CH}^{(2)}_{s}\left(
\mathfrak{Q}'^{(2)}
\right)
\\&=
\mathrm{CH}^{(2)}_{s}\left(
\mathfrak{P}^{(2),i+1,\bb{\mathfrak{P}^{(2)}}-1}
\right)
\circ^{1\mathbf{T}_{\Sigma^{\boldsymbol{\mathcal{A}}^{(2)}}}(X)}_{s}
\mathrm{CH}^{(2)}_{s}\left(
\mathfrak{Q}^{(2)}
\right)
\\&=
\mathrm{CH}^{(2)}_{s}\left(\mathfrak{P}^{(2)}
\right).
\end{align*}

The case $\mathfrak{Q}^{(2)}=\mathfrak{P}^{(2),0,i}$ follows.

If~(3.2.2), i.e., if $\mathfrak{P}^{(2)}$ is a head-constant echelonless second-order path that is not coherent, $i\in\bb{\mathfrak{P}^{(2)}}$ is the maximum index for which $\mathfrak{P}^{(2),0,i}$ is a coherent head-constant echelonless second-order path and $\mathfrak{Q}^{(2)}=\mathfrak{P}^{(2),i+1,\bb{\mathfrak{P}^{(2)}}-1}$ then, following Definition~\ref{DDCH}, the value of the second-order Curry-Howard mapping at $\mathfrak{P}^{(2)}$ is given by
$$
\mathrm{CH}^{(2)}_{s}\left(
\mathfrak{P}^{(2)}
\right)
=
\mathrm{CH}^{(2)}_{s}\left(
\mathfrak{Q}^{(2)}
\right)
\circ^{1\mathbf{T}_{\Sigma^{\boldsymbol{\mathcal{A}}^{(2)}}}(X)}_{s}
\mathrm{CH}^{(2)}_{s}\left(
\mathfrak{P}^{(2),0,i}
\right).
$$

Since $\mathfrak{Q}^{(2)}$ is a head-constant echelonless second-order path of length, then in virtue of Lemma~\ref{LDCHNEchHd} we have that $\mathrm{CH}^{(2)}_{s}(\mathfrak{Q}^{(2)})\in [\mathrm{T}_{\Sigma^{\boldsymbol{\mathcal{A}}^{(2)}}}(X)]^{\mathsf{HdC}}_{s}
.$ 
Thus, if $\mathfrak{Q}'^{(2)}$ is any-second-order path in $[\mathfrak{Q}^{(2)}]^{}_{s}$ then, in virtue of Lemma~\ref{LDCHNEchHd} we conclude that $\mathfrak{Q}'^{(2)}$  is a head-constant echelonless second-order path. Moreover, in virtue of Lemma~\ref{LDCH}, we have that 
$
\mathrm{sc}^{([1],2)}_{s}(\mathfrak{Q}'^{(2)})
=
\mathrm{sc}^{([1],2)}_{s}(\mathfrak{Q}^{(2)})
$ and the length of $\mathfrak{Q}'^{(2)}$ equals the length of $\mathfrak{Q}^{(2)}$. Therefore, since $\mathfrak{P}^{(2)}$ decomposes as the $1$-composition $\mathfrak{P}^{(2)}=\mathfrak{Q}^{(2)}\circ^{1\mathbf{Pth}_{\boldsymbol{\mathcal{A}}^{(2)}}}_{s}\mathfrak{P}^{(2),0,i}$, we have that 
$
\mathrm{sc}^{([1],2)}_{s}(\mathfrak{Q}^{(2)})
=
\mathrm{tg}^{([1],2)}_{s}(\mathfrak{P}^{(2),0,i}).
$ All in all, we conclude that 
$
\mathrm{sc}^{([1],2)}_{s}(\mathfrak{Q}'^{(2)})
=
\mathrm{tg}^{([1],2)}_{s}(\mathfrak{P}^{(2),0,i}).
$ 
Therefore the second-order path $\mathfrak{P}'^{(2)}$ defined by 
$$
\mathfrak{P}'^{(2)}
=\mathfrak{Q}'^{(2)}\circ^{1\mathbf{Pth}_{\boldsymbol{\mathcal{A}}^{(2)}}}_{s}\mathfrak{P}^{(2),0,i}
$$
is, in virtue of Lemma~\ref{LTech}, a head-constant echelonless second-order path that is not coherent for which $i\in\bb{\mathfrak{P}^{(2)}}$ is the maximum index for which $\mathfrak{P}'^{(2),0,i}=\mathfrak{P}^{(2),0,i}$ is a coherent head-constant echelonless second-order path.  Taking into account Definition~\ref{DDOrd}, we have that
$
(\mathfrak{Q}'^{(2)},s)
\prec_{\mathbf{Pth}_{\boldsymbol{\mathcal{A}}^{(2)}}}
(\mathfrak{P}'^{(2)},s).
$ Moreover the value of the second-order Curry-Howard mapping at $\mathfrak{P}'^{(2)}$ is given by
\begin{align*}
\mathrm{CH}^{(2)}_{s}\left(
\mathfrak{P}'^{(2)}
\right)
&=
\mathrm{CH}^{(2)}_{s}\left(
\mathfrak{Q}'^{(2)}
\right)
\circ^{1\mathbf{T}_{\Sigma^{\boldsymbol{\mathcal{A}}^{(2)}}}(X)}_{s}
\mathrm{CH}^{(2)}_{s}\left(
\mathfrak{P}^{(2),0,i}
\right)
\\&=
\mathrm{CH}^{(2)}_{s}\left(
\mathfrak{Q}^{(2)}
\right)
\circ^{1\mathbf{T}_{\Sigma^{\boldsymbol{\mathcal{A}}^{(2)}}}(X)}_{s}
\mathrm{CH}^{(2)}_{s}\left(
\mathfrak{P}^{(2),0,i}
\right)
\\&=
\mathrm{CH}^{(2)}_{s}\left(
\mathfrak{P}^{(2)}
\right).
\end{align*}

The case $\mathfrak{Q}^{(2)}=\mathfrak{P}^{(2),i+1,\bb{\mathfrak{P}^{(2)}}-1}$ follows.

This completes the case (3.2) of $\mathfrak{P}^{(2)}$ being a head-constant echelonless second-order path  that is not coherent. 

Finally, consider the case (3.3), in which $\mathfrak{P}^{(2)}$ is a coherent head-constant echelonless second-order path. Then, following Definition~\ref{DDHeadCt}, there exists a  unique word $\mathbf{s}\in S^{\star}-\{\lambda\}$ and a unique operation symbol $\tau\in\Sigma^{\boldsymbol{\mathcal{A}}}_{\mathbf{s},s}$ associated to $\mathfrak{P}^{(2)}$. Let $(\mathfrak{P}^{(2)}_{j})_{j\in\bb{\mathbf{s}}}$ be the family of second-order paths we can extract from $\mathfrak{P}^{(2)}$ in virtue of Lemma~\ref{LDPthExtract}. Assume that, for an index $j\in\bb{\mathbf{s}}$, it is the case that $\mathfrak{Q}^{(2)}=\mathfrak{P}^{(2)}_{j}$. Then, following Definition~\ref{DDCH}, the value of the second-order Curry-Howard mapping at $\mathfrak{P}^{(2)}$ is given by
$$
\mathrm{CH}^{(2)}_{s}
\left(
\mathfrak{P}^{(2)}
\right)
=
\tau^{\mathbf{T}_{\Sigma^{\boldsymbol{\mathcal{A}}^{(2)}}}(X)}
\left(\left(\mathrm{CH}^{(2)}_{s_{j}}\left(
\mathfrak{P}^{(2)}_{j}
\right)\right)_{j\in\bb{\mathbf{s}}}
\right).
$$

Thus, if $\mathfrak{Q}'^{(2)}$ is any second-order path in $[\mathfrak{Q}^{(2)}]^{}_{s_{j}}$ then, in virtue of Lemma~\ref{LDCH}, we have that $\mathfrak{Q}'^{(2)}$ and $\mathfrak{Q}^{(2)}$ have the same length, $([1],2)$-source and $([1],2)$-target.

Consider the second-order path $\mathfrak{P}'^{(2)}$ in $\mathrm{Pth}_{\boldsymbol{\mathcal{A}}^{(2)},s}$, given by
$$
\mathfrak{P}'^{(2)}
=
\tau^{\mathbf{Pth}_{\boldsymbol{\mathcal{A}}^{(2)}}}
\left(\mathfrak{P}^{(2)}_{0},
\cdots,
\mathfrak{Q}'^{(2)},
\cdots,
\mathfrak{P}^{(2)}_{\bb{\mathbf{s}}-1}
\right),
$$
that is, $\mathfrak{P}'^{(2)}$ is defined as the operation $\tau^{\mathbf{Pth}_{\boldsymbol{\mathcal{A}}^{(2)}}}$ applied to the family of second-order paths extracted from $\mathfrak{P}^{(2)}$, where we have replaced  the second-order path at position $j$ by $\mathfrak{Q}'^{(2)}$. 

Note that, for the special case of $\tau$ being the operation symbol of $0$-composition $\circ^{0}_{s}$, the second-order path $\mathfrak{P}'^{(2)}$ is well-defined, because $\mathfrak{Q}'^{(2)}$ and $\mathfrak{Q}^{(2)}$ have the same $([1],2)$-source and $([1],2)$-target. In particular, following Definition~\ref{DDScTgZ}, they have the same $(0,2)$-source and $(0,2)$-target. Thus, they are mutually interchangeable in any $0$-composition.

Note that, in virtue of Proposition~\ref{PDRecov}, $\mathfrak{P}'^{(2)}$ is a  
coherent head-constant echelonless second-order path for which, the second-order path extraction procedure applied to $\mathfrak{P}'^{(2)}$ retrieves the family of second-order paths $(\mathfrak{P}^{(2)},\cdots, \mathfrak{Q}'^{(2)},\cdots, \mathfrak{P}^{(2)}_{\bb{\mathbf{s}}-1})$. Therefore, following Definition~\ref{DDOrd}, we conclude that 
$(\mathfrak{Q}'^{(2)},s_{j})
\prec_{\mathbf{Pth}_{\boldsymbol{\mathcal{A}}^{(2)}}}
(\mathfrak{P}'^{(2)},s).
$
Moreover, the value of the second-order Curry-Howard mapping at $\mathfrak{P}^{(2)}$ is given by 
\begin{flushleft}
$\mathrm{CH}^{(2)}_{s}\Bigr(
\mathfrak{P}'^{(2)}
\Bigr)$
\allowdisplaybreaks
\begin{align*}
\qquad
&=
\tau^{\mathbf{T}_{\Sigma^{\boldsymbol{\mathcal{A}}^{(2)}}}(X)}
\left(\mathrm{CH}^{(2)}_{s_{0}}\left(
\mathfrak{P}^{(2)}_{0}
\right),
\cdots,
\mathrm{CH}^{(2)}_{s_{j}}\left(
\mathfrak{Q}'^{(2)}
\right),
\cdots,
\mathrm{CH}^{(2)}_{s_{\bb{\mathbf{s}}-1}}\left(
\mathfrak{P}^{(2)}_{\bb{\mathbf{s}}-1}
\right)\right)
\\&=
\tau^{\mathbf{T}_{\Sigma^{\boldsymbol{\mathcal{A}}^{(2)}}}(X)}
\left(\left(\mathrm{CH}^{(2)}_{s_{j}}\left(
\mathfrak{P}^{}_{j}
\right)\right)_{j\in\bb{\mathbf{s}}}
\right)
\\&=
\mathrm{CH}^{(2)}_{s}\left(
\mathfrak{P}^{(2)}
\right).
\end{align*}
\end{flushleft}

This completes the case (3.3) of $\mathfrak{P}^{(2)}$ being a coherent head-constant echelonless second-order path.

This completes the Case~(3).

This concludes the base case.

\textsf{Inductive step of the induction.}

Assume the statement holds for sequences of length up to $m$, i.e., for every pair of sorts $t,s\in S$, if $\mathfrak{Q}'^{(2)},\mathfrak{Q}^{(2)}$ are second-order paths in $\mathrm{Pth}_{\boldsymbol{\mathcal{A}}^{(2)},t}$ and $\mathfrak{P}^{(2)}$ is a second-order path in $\mathrm{Pth}_{\boldsymbol{\mathcal{A}}^{(2)},s}$ such that  there exists a word $\mathbf{w}\in S^{\star}$ of length $\bb{\mathbf{w}}=m+1$ and a family of second-order paths $(\mathfrak{R}^{(2)}_{k})_{k\in\bb{\mathbf{w}}}$ in $\mathrm{Pth}_{\boldsymbol{\mathcal{A}}^{(2)},\mathbf{w}}$  such that $w_{0}=t$, $\mathfrak{R}^{(2)}_{0}=\mathfrak{Q}^{(2)}$, $w_{m}=s$, $\mathfrak{R}^{(2)}_{m}=\mathfrak{P}^{(2)}$ and 
for every $k\in m$, $(\mathfrak{R}^{(2)}_{k}, w_{k})\prec_{\mathbf{Pth}_{\boldsymbol{\mathcal{A}}^{(2)}}} (\mathfrak{R}^{(2)}_{k+1}, w_{k+1})$ and  $\mathfrak{Q}'^{(2)}$ is a second-order path in $[\mathfrak{Q}^{(2)}]^{}_{t}$, then there exists a second-order path $\mathfrak{P}'^{(2)}$ in $[\mathfrak{P}^{(2)}]^{}_{s}$ such that 
$$
\left(
\mathfrak{Q}'^{(2)},t
\right)\leq_{\mathbf{Pth}_{\boldsymbol{\mathcal{A}}^{(2)}}}
\left(
\mathfrak{P}'^{(2)},s
\right).
$$

Now we prove it for sequences of length $m+1$. 

Let $t,s$ be sorts in $S$, if $\mathfrak{Q}'^{(2)},\mathfrak{Q}^{(2)}$ are second-order paths in $\mathrm{Pth}_{\boldsymbol{\mathcal{A}}^{(2)},t}$ and $\mathfrak{P}^{(2)}$ is a second-order path in $\mathrm{Pth}_{\boldsymbol{\mathcal{A}}^{(2)},s}$ such that  there exists a word $\mathbf{w}\in S^{\star}$ of length $\bb{\mathbf{w}}=m+2$ and a family of second-order paths $(\mathfrak{R}^{(2)}_{k})_{k\in\bb{\mathbf{w}}}$ in $\mathrm{Pth}_{\boldsymbol{\mathcal{A}}^{(2)},\mathbf{w}}$  such that $w_{0}=t$, $\mathfrak{R}^{(2)}_{0}=\mathfrak{Q}^{(2)}$, $w_{m+1}=s$, $\mathfrak{R}^{(2)}_{m+1}=\mathfrak{P}^{(2)}$ and 
for every $k\in m+1$, $(\mathfrak{R}^{(2)}_{k}, w_{k})\prec_{\mathbf{Pth}_{\boldsymbol{\mathcal{A}}^{(2)}}} (\mathfrak{R}^{(2)}_{k+1}, w_{k+1})$ and let $\mathfrak{Q}'^{(2)}$ be a second-order path in $[\mathfrak{Q}^{(2)}]^{}_{t}$.

Consider the word $\mathbf{w}^{0,m}$ of length $\bb{\mathbf{w}^{0,m}}=m+1$ and the family of second-order paths $(\mathfrak{R}^{(2)}_{k})_{k\in\bb{\mathbf{w}^{0,m}}}$ in $\mathrm{Pth}_{\boldsymbol{\mathcal{A}}^{(2)},\mathbf{w}^{0,m}}$. This is a sequence of length $m$ instantiating that 
$
(\mathfrak{Q}^{(2)},t)
\leq_{\mathbf{Pth}_{\boldsymbol{\mathcal{A}}^{(2)}}}
(\mathfrak{R}^{(2)}_{m}, w_{m}).
$ Since $\mathfrak{Q}'^{(2)}$ is a second-order path in $[\mathfrak{Q}^{(2)}]^{}_{t}$, by the inductive hypothesis, there exist a second-order path $\mathfrak{R}'^{(2)}_{m}$ in $[\mathfrak{R}^{(2)}_{m}]^{}_{w_{m}}$ satisfying that 
$
(\mathfrak{Q}'^{(2)},t)
\leq_{\mathbf{Pth}_{\boldsymbol{\mathcal{A}}^{(2)}}}
(\mathfrak{R}'^{(2)}_{m}, w_{m}).
$

Now, consider the sequence of second-order paths $(\mathfrak{R}^{(2)}_{m},\mathfrak{P}^{(2)})$. It is a one-step sequence of second-order paths in $\mathrm{Pth}_{\boldsymbol{\mathcal{A}}^{(2)},\mathbf{w}^{m,m+1}}$ satisfying that $(\mathfrak{R}^{(2)}_{m},w_{m})\prec_{\mathbf{Pth}_{\boldsymbol{\mathcal{A}}^{(2)}}} (\mathfrak{P}^{(2)}, s)$. Since $\mathfrak{R}'^{(2)}_{m}$ in $[\mathfrak{R}^{(2)}_{m}]^{}_{w_{m}}$, by the base case, we can find a second-order path $\mathfrak{P}'^{(2)}$ in $[\mathfrak{P}^{(2)}]^{}_{s}$ satisfying that 
$
(\mathfrak{R}'^{(2)}_{m},w_{m})
\leq_{\mathbf{Pth}_{\boldsymbol{\mathcal{A}}^{(2)}}}
(\mathfrak{P}'^{(2)}, s).
$

Hence, $\mathfrak{P}'^{(2)}$ is a second-order path in $[\mathfrak{P}^{(2)}]^{}_{s}$ satisfying that 
$$
\left(
\mathfrak{Q}'^{(2)},t
\right)
\leq_{\mathbf{Pth}_{\boldsymbol{\mathcal{A}}^{(2)}}}
\left(
\mathfrak{P}'^{(2)}, s
\right).
$$

This concludes the proof.
\end{proof}

\begin{corollary}\label{CDCHOrdI} Let $t,s$ be sorts in $S$, $[\mathfrak{Q}^{(2)}]^{}_{t}$ a second-order path class in $[\mathrm{Pth}_{\boldsymbol{\mathcal{A}}^{(2)}}]_{t}$ and  $[\mathfrak{P}^{(2)}]^{}_{s}$ a second-order path class in $[\mathrm{Pth}_{\boldsymbol{\mathcal{A}}^{(2)}}]_{s}$. Then the following statements are equivalent
\begin{enumerate}
\item[(i)] $([\mathfrak{Q}^{(2)}]^{}_{t},t)
\leq_{[\mathbf{Pth}_{\boldsymbol{\mathcal{A}}^{(2)}}]}
([\mathfrak{P}^{(2)}]^{}_{s},s);
$
\item[(ii)] There exists $
\mathfrak{P}'^{(2)}\in [\mathfrak{P}^{(2)}]^{}_{s}$ such that $
(\mathfrak{Q}^{(2)},t)
\leq_{\mathbf{Pth}_{\boldsymbol{\mathcal{A}}^{(2)}}}
(\mathfrak{P}'^{(2)},s)
.
$
\end{enumerate}
\end{corollary}
\begin{proof}
Assume that $([\mathfrak{Q}^{(2)}]^{}_{t},t)
\leq_{[\mathbf{Pth}_{\boldsymbol{\mathcal{A}}^{(2)}}]}
([\mathfrak{P}^{(2)}]^{}_{s},s)
$ then there exists $\mathfrak{Q}'^{(2)}\in[\mathfrak{Q}^{(2)}]^{}_{t}$ and $\mathfrak{P}'^{(2)}\in [\mathfrak{P}^{(2)}]^{}_{s}$ such that $(\mathfrak{Q}'^{(2)},t) \leq_{\mathbf{Pth}_{\boldsymbol{\mathcal{A}}^{(2)}}} (\mathfrak{P}'^{(2)},s)$. Since $\mathfrak{Q}^{(2)}\in[\mathfrak{Q}'^{(2)}]^{}_{t}$, by Lemma~\ref{LDCHOrdI}, we can find $\mathfrak{P}''^{(2)}\in[\mathfrak{P}'^{(2)}]^{}_{s}$ such that $(\mathfrak{Q}^{(2)},t) \leq_{\mathbf{Pth}_{\boldsymbol{\mathcal{A}}^{(2)}}} (\mathfrak{P}''^{(2)},s)$. Note that $\mathfrak{P}''^{(2)}\in[\mathfrak{P}^{(2)}]^{}_{s}$.

The other implication follows by definition of the relation $\leq_{[\mathbf{Pth}_{\boldsymbol{\mathcal{A}}^{(2)}}]}$.

This concludes the proof.
\end{proof}

Next lemma shows that, if we are given two second-order paths in the same class under the Kernel of the second-order Curry Howard mapping and one of these second-order paths is greater than another second-order path with respect to the partial order $\leq_{\mathbf{Pth}_{\boldsymbol{\mathcal{A}}^{(2)}}}$, then we can complete this inequality on the other side with an equivalent second-order path.

\begin{lemma}\label{LDCHOrdII} Let $t,s$ be sorts in $S$, $\mathfrak{Q}^{(2)}$ a second-order path in $\mathrm{Pth}_{\boldsymbol{\mathcal{A}}^{(2)},t}$ and $\mathfrak{P}'^{(2)},\mathfrak{P}^{(2)}$ second-order paths in $\mathrm{Pth}_{\boldsymbol{\mathcal{A}}^{(2)},s}$. If  $(\mathfrak{Q}^{(2)},t)\leq_{\mathbf{Pth}_{\boldsymbol{\mathcal{A}}^{(2)}}}(\mathfrak{P}^{(2)},s)$ and $\mathfrak{P}'^{(2)}\in[\mathfrak{P}^{(2)}]^{}_{s}$ then there exists a second-order path $\mathfrak{Q}'^{(2)}\in [\mathfrak{Q}^{(2)}]^{}_{t}$ satisfying that  $$
\left(
\mathfrak{Q}'^{(2)},t
\right)\leq_{\mathbf{Pth}_{\boldsymbol{\mathcal{A}}^{(2)}}}
\left(
\mathfrak{P}'^{(2)},s
\right).$$
\end{lemma}

\begin{proof}
Let us recall from Remark~\ref{RDOrd} that $(\mathfrak{Q}^{(2)},t)\leq_{\mathbf{Pth}_{\boldsymbol{\mathcal{A}}^{(2)}}}(\mathfrak{P}^{(2)},s)$ if and only if $s=t$ and $\mathfrak{Q}^{(2)}=\mathfrak{P}^{(2)}$ or there exists a natural number $m\in\mathbb{N}-\{0\}$, a word $\mathbf{w}\in S^{\star}$ of length $\bb{\mathbf{w}}=m+1$, and a family of second-order paths $(\mathfrak{R}^{(2)}_{k})$ in $\mathrm{Pth}_{\boldsymbol{\mathcal{A}}^{(2)},\mathbf{w}}$ such that $w_{0}=t$, $\mathfrak{R}^{(2)}_{0}=\mathfrak{Q}^{(2)}$, $w_{m}=s$, $\mathfrak{R}^{(2)}_{m}=\mathfrak{P}^{(2)}$ and for every $k\in m$, $(\mathfrak{R}^{(2)}_{k}, w_{k})\prec_{\mathbf{Pth}_{\boldsymbol{\mathcal{A}}^{(2)}}} (\mathfrak{R}^{(2)}_{k+1}, w_{k+1})$.

The lemma holds trivially in case $s=t$ and $\mathfrak{Q}^{(2)}=\mathfrak{P}^{(2)}$. Therefore, it remains to prove the other case. 

We will prove the lemma by induction on $m\in\mathbb{N}-\{0\}$.

\textsf{Base step of the induction.}

For $m=1$ we have that $(\mathfrak{Q}^{(2)},t)\prec_{\mathbf{Pth}_{\boldsymbol{\mathcal{A}}^{(2)}}}(\mathfrak{P}^{(2)},s)$, hence following Definition~\ref{DDOrd} we are in one of the following situations
\begin{enumerate}
\item $\mathfrak{P}^{(2)}$ and $\mathfrak{Q}^{(2)}$ are $(2,[1])$-identity second-order paths of the form 
\begin{align*}
\mathfrak{P}^{(2)}&=
\mathrm{ip}^{(2,[1])\sharp}_{s}\left(
[P]_{s}
\right),
&
\mathfrak{Q}^{(2)}&=
\mathrm{ip}^{(2,[1])\sharp}_{t}\left(
[Q]_{t}
\right),
\end{align*}
for some path term classes $[P]_{s}\in[\mathrm{PT}_{\boldsymbol{\mathcal{A}}}]_{s}$ and  $[Q]_{t}\in[\mathrm{PT}_{\boldsymbol{\mathcal{A}}}]_{t}$ and the following inequality holds
$$
\left(
[Q]_{t},t
\right)
<_{[\mathbf{PT}_{\boldsymbol{\mathcal{A}}}]
}
\left(
[P]_{s},s
\right),
$$
where $\leq_{[\mathbf{PT}_{\boldsymbol{\mathcal{A}}}]
}$ is the Artinian partial order on $\coprod [\mathrm{PT}_{\boldsymbol{\mathcal{A}}}]$ introduced in Definition~\ref{DPTQOrd}.
\item $\mathfrak{P}^{(2)}$ is a second-order path of length strictly greater than one containing at least one   second-order echelon, and if its first   second-order echelon occurs at position $i\in\bb{\mathfrak{P}^{(2)}}$, then
\subitem if $i=0$, then $\mathfrak{Q}^{(2)}$ is equal to $\mathfrak{P}^{(2),0,0}$ or $\mathfrak{P}^{(2),1,\bb{\mathfrak{P}^{(2)}}-1}$,
\subitem if $i>0$, then $\mathfrak{Q}^{(2)}$ is equal to $\mathfrak{P}^{(2),0,i-1}$ or $\mathfrak{P}^{(2),i,\bb{\mathfrak{P}^{(2)}}-1}$,
\item $\mathfrak{P}^{(2)}$ is an echelonless second-order path of length at least one, then
\subitem if $\mathfrak{P}^{(2)}$ is not head-constant, then let $i\in\bb{\mathfrak{P}^{(2)}}$ be the  maximum index for which $\mathfrak{P}^{(2),0,i}$ is a head-constant second-order path, then $\mathfrak{Q}^{(2)}$ is equal to $\mathfrak{P}^{(2),0,i}$ or $\mathfrak{P}^{(2),i+1,\bb{\mathfrak{P}^{(2)}}-1}$,
\subitem if $\mathfrak{P}^{(2)}$ is head-constant but not coherent, then let $i\in\bb{\mathfrak{P}^{(2)}}$ be the  maximum index for which $\mathfrak{P}^{(2),0,i}$ is a coherent second-order path, then $\mathfrak{Q}^{(2)}$ is equal to $\mathfrak{P}^{(2),0,i}$ or $\mathfrak{P}^{(2),i+1,\bb{\mathfrak{P}^{(2)}}-1}$,
\subitem if $\mathfrak{P}^{(2)}$ is head-constant and coherent then $\mathfrak{Q}^{(2)}$ is one of the second-order paths we can extract from $\mathfrak{P}^{(2)}$ in virtue of Lemma~\ref{LDPthExtract}.
\end{enumerate}

If~(1), we note that both $\mathfrak{P}^{(2)}$ and $\mathfrak{Q}^{(2)}$ are $(2,[1])$-identity second-order paths. Then, if $\mathfrak{P}'^{(2)}$ is a second-order path in $[\mathfrak{P}^{(2)}]^{}_{s}$ then we conclude, in virtue  of Corollary~\ref{CDCHUId}, that $\mathfrak{P}'^{(2)}=\mathfrak{P}^{(2)}$. This case follows easily, since $\mathfrak{Q}^{(2)}$ is a second-order path in $[\mathfrak{Q}^{(2)}]^{}_{t}$ satisfying that  $$ \left(\mathfrak{Q}^{(2)},t
\right)\leq_{\mathbf{Pth}_{\boldsymbol{\mathcal{A}}^{(2)}}}
\left(
\mathfrak{P}^{(2)},s
\right).$$

This completes the Case~(1).

If~(2), that is, if $\mathfrak{P}^{(2)}$ is a second-order path of length strictly greater than one containing at least one   second-order echelon and its first   second-order echelon occurs at position $i\in\bb{\mathfrak{P}^{(2)}}$, we consider the case~(2.1) $i=0$ and the case~(2.2) $i>0$ .

If~(2.1), i.e.,  if $\mathfrak{P}^{(2)}$ is a second-order path of length strictly greater than one containing its first   second-order echelon on its first step then, following Definition~\ref{DDCH}, the value of the second-order Curry-Howard mapping at $\mathfrak{P}^{(2)}$ is given by
$$
\mathrm{CH}^{(2)}_{s}\left(
\mathfrak{P}^{(2)}
\right)
=
\mathrm{CH}^{(2)}_{s}\left(
\mathfrak{P}^{(2),1,\bb{\mathfrak{P}^{(2)}}-1}
\right)
\circ^{1\mathbf{T}_{\Sigma^{\boldsymbol{\mathcal{A}}^{(2)}}}(X)}_{s}
\mathrm{CH}^{(2)}_{s}\left(
\mathfrak{P}^{(2),0,0}
\right).
$$

Since $\mathfrak{P}^{(2)}$ is a second-order path of length strictly greater than one containing a   second-order echelon on its first step then, in virtue of Lemma~\ref{LDCHEchInt}, we have that 
$\mathrm{CH}^{(2)}_{s}(\mathfrak{P}^{(2)})\in \eta^{(2,\mathcal{A}^{(2)})}[\mathcal{A}^{(2)}]^{\mathrm{int}}_{s}$. Thus, if $\mathfrak{P}'^{(2)}$ is any second-order path in $[\mathfrak{P}^{(2)}]^{}_{s}$ then,  in virtue of Lemma~\ref{LDCHEchInt}, we have that  $\mathfrak{P}'^{(2)}$ is a second-order path of length strictly greater than one containing a   second-order echelon on its first step.

Thus, the value of the second-order Curry-Howard mapping at $\mathfrak{P}'^{(2)}$ is given by
$$
\mathrm{CH}^{(2)}_{s}\left(
\mathfrak{P}'^{(2)}
\right)
=
\mathrm{CH}^{(2)}_{s}\left(
\mathfrak{P}'^{(2),1,\bb{\mathfrak{P}'^{(2)}}-1}
\right)
\circ^{1\mathbf{T}_{\Sigma^{\boldsymbol{\mathcal{A}}^{(2)}}}(X)}_{s}
\mathrm{CH}^{(2)}_{s}\left(
\mathfrak{P}'^{(2),0,0}
\right).
$$

Now, we consider the different possibilities for $\mathfrak{Q}^{(2)}$. Note that, since $(\mathfrak{Q}^{(2)}, t)\prec_{\mathbf{Pth}_{\boldsymbol{\mathcal{A}}^{(2)}}} (\mathfrak{P}^{(2)},s)$, then following Definition~\ref{DDOrd}, either (2.1.1) $\mathfrak{Q}^{(2)}=\mathfrak{P}^{(2),0,0}$ or~(2.1.2) $\mathfrak{Q}^{(2)}=\mathfrak{P}^{(2),1,\bb{\mathfrak{P}^{(2)}}-1}$. In any case, note that $t=s$.

If~(2.1.1), i.e., if $\mathfrak{Q}^{(2)}=\mathfrak{P}^{(2),0,0}$ then, following Definition~\ref{DDOrd}, $(\mathfrak{P}'^{(2),0,0},s)\prec_{\mathbf{Pth}_{\boldsymbol{\mathcal{A}}^{(2)}}}(\mathfrak{P}'^{(2)},s)$. Moreover, since $\mathfrak{P}'^{(2)}$ is a second-order path in $[\mathfrak{P}^{(2)}]^{}_{s}$ we conclude that $\mathrm{CH}^{(2)}_{s}(\mathfrak{Q}^{(2)})=\mathrm{CH}^{(2)}_{s}(\mathfrak{P}'^{(2),0,0})$.

On the other hand, if~(2.1.2), i.e., if $\mathfrak{Q}^{(2)}=\mathfrak{P}^{(2),1,\bb{\mathfrak{P}^{(2)}}-1}$ then, following Definition~\ref{DDOrd}, $(\mathfrak{P}'^{(2),1,\bb{\mathfrak{P}'^{(2)}}-1},s)\prec_{\mathbf{Pth}_{\boldsymbol{\mathcal{A}}^{(2)}}}(\mathfrak{P}'^{(2)},s)$. Moreover, since $\mathfrak{P}'^{(2)}$ is a second-order path in $[\mathfrak{P}^{(2)}]^{}_{s}$ we conclude that $\mathrm{CH}^{(2)}_{s}(\mathfrak{Q}^{(2)})=\mathrm{CH}^{(2)}_{s}(\mathfrak{P}'^{(2),1,\bb{\mathfrak{P}^{(2)}}-1})$.

The case $i=0$ follows.

If~(2.2), i.e.,  if $\mathfrak{P}^{(2)}$ is a second-order path of length strictly greater than one containing its first   second-order echelon on the step $i\in\bb{\mathfrak{P}^{(2)}}$ and $i>0$, then following Definition~\ref{DDCH}, the value of the second-order Curry-Howard mapping at $\mathfrak{P}^{(2)}$ is given by
$$
\mathrm{CH}^{(2)}_{s}\left(
\mathfrak{P}^{(2)}
\right)
=
\mathrm{CH}^{(2)}_{s}\left(
\mathfrak{P}^{(2),i,\bb{\mathfrak{P}^{(2)}}-1}
\right)
\circ^{1\mathbf{T}_{\Sigma^{\boldsymbol{\mathcal{A}}^{(2)}}}(X)}_{s}
\mathrm{CH}^{(2)}_{s}\left(
\mathfrak{P}^{(2),0,i-1}
\right).
$$

Since $\mathfrak{P}^{(2)}$ is a second-order path of length strictly greater than one containing a   second-order echelon on a step different from the initial one then, in virtue of Lemma~\ref{LDCHEchNInt}, we have that 
$\mathrm{CH}^{(2)}_{s}(\mathfrak{P}^{(2)})\in \eta^{(2,\mathcal{A}^{(2)})}[\mathcal{A}^{(2)}]^{\neg\mathrm{int}}_{s}$. Thus, if $\mathfrak{P}'^{(2)}$ is any second-order path in $[\mathfrak{P}^{(2)}]^{}_{s}$ then,  in virtue of Lemma~\ref{LDCHEchNInt}, we have that  $\mathfrak{P}'^{(2)}$ is a second-order path of length strictly greater than one containing a   second-order echelon on a step different from the initial one. Let $j\in\bb{\mathfrak{P}'^{(2)}}$ be the first index for which $\mathfrak{P}'^{(2),j,j}$ is a   second-order echelon.

Thus, the value of the second-order Curry-Howard mapping at $\mathfrak{P}'^{(2)}$ is given by
$$
\mathrm{CH}^{(2)}_{s}\left(
\mathfrak{P}'^{(2)}
\right)
=
\mathrm{CH}^{(2)}_{s}\left(
\mathfrak{P}'^{(2),j,\bb{\mathfrak{P}'^{(2)}}-1}
\right)
\circ^{1\mathbf{T}_{\Sigma^{\boldsymbol{\mathcal{A}}^{(2)}}}(X)}_{s}
\mathrm{CH}^{(2)}_{s}\left(
\mathfrak{P}'^{(2),0,j-1}
\right).
$$

Now, we consider the different possibilities for $\mathfrak{Q}^{(2)}$. Note that, since $(\mathfrak{Q}^{(2)}, t)\prec_{\mathbf{Pth}_{\boldsymbol{\mathcal{A}}^{(2)}}} (\mathfrak{P}^{(2)},s)$, then following Definition~\ref{DDOrd}, either~(2.2.1) $\mathfrak{Q}^{(2)}=\mathfrak{P}^{(2),0,i-1}$ or~(2.2.2) $\mathfrak{Q}^{(2)}=\mathfrak{P}^{(2),i,\bb{\mathfrak{P}^{(2)}}-1}$. In any case, note that $t=s$.

If~(2.2.1), i.e., if $\mathfrak{Q}^{(2)}=\mathfrak{P}^{(2),0,i-1}$ then, following Definition~\ref{DDOrd}, $(\mathfrak{P}'^{(2),0,j-1},s)\prec_{\mathbf{Pth}_{\boldsymbol{\mathcal{A}}^{(2)}}}(\mathfrak{P}'^{(2)},s)$. Moreover, since $\mathfrak{P}'^{(2)}$ is a second-order path in $[\mathfrak{P}^{(2)}]^{}_{s}$ we conclude that $\mathrm{CH}^{(2)}_{s}(\mathfrak{Q}^{(2)})=\mathrm{CH}^{(2)}_{s}(\mathfrak{P}'^{(2),0,j-1})$.

On the other hand, if~(2.2.2), i.e., if $\mathfrak{Q}^{(2)}=\mathfrak{P}^{(2),i,\bb{\mathfrak{P}^{(2)}}-1}$ then, following Definition~\ref{DDOrd}, $(\mathfrak{P}'^{(2),j,\bb{\mathfrak{P}'^{(2)}}-1},s)\prec_{\mathbf{Pth}_{\boldsymbol{\mathcal{A}}^{(2)}}}(\mathfrak{P}'^{(2)},s)$. Moreover, since $\mathfrak{P}'^{(2)}$ is a second-order path in $[\mathfrak{P}^{(2)}]^{}_{s}$ we conclude that $\mathrm{CH}^{(2)}_{s}(\mathfrak{Q}^{(2)})=\mathrm{CH}^{(2)}_{s}(\mathfrak{P}'^{(2),j,\bb{\mathfrak{P}^{(2)}}-1})$.

The case $i>0$ follows.

This completes the Case~(2).

If~(3), that is, if $\mathfrak{P}^{(2)}$ is an echelonless second-order path, then we consider different cases according to the nature of $\mathfrak{P}^{(2)}$. In this regard, we consider the case (3.1) $\mathfrak{P}^{(2)}$ is an echelonless second-order path that is not head-constant, (3.2) $\mathfrak{P}^{(2)}$ is a head-constant echelonless second-order path that is not coherent, and~(3.3) $\mathfrak{P}^{(2)}$ is a coherent head-constant echelonless second-order path.

If~(3.1), i.e., if $\mathfrak{P}^{(2)}$ is an echelonless second-order path that is not head-constant, then let $i\in\bb{\mathfrak{P}^{(2)}}$ be the maximum index for which $\mathfrak{P}^{(2),0,i}$ is a head-constant echelonless second-order path. Following Definition~\ref{DDCH}, the value of the second-order Curry-Howard mapping at $\mathfrak{P}^{(2)}$ is given by
$$
\mathrm{CH}^{(2)}_{s}\left(
\mathfrak{P}^{(2)}
\right)
=
\mathrm{CH}^{(2)}_{s}\left(
\mathfrak{P}^{(2),i+1,\bb{\mathfrak{P}^{(2)}}-1}
\right)
\circ^{1\mathbf{T}_{\Sigma^{\boldsymbol{\mathcal{A}}^{(2)}}}(X)}_{s}
\mathrm{CH}^{(2)}_{s}\left(
\mathfrak{P}^{(2),0,i}
\right).
$$

Since $\mathfrak{P}^{(2)}$ is an echelonless second-order path of length at least one that is not head-constant then, in virtue of Lemma~\ref{LDCHNEchNHd}, we have that $\mathrm{CH}^{(2)}_{s}(\mathfrak{P}^{(2)})\in 
[\mathrm{T}_{\Sigma^{\boldsymbol{\mathcal{A}}^{(2)}}}(X)]^{\neg\mathsf{HdC}}_{s}
$.  Thus, if $\mathfrak{P}'^{(2)}$ is any second-order path in $[\mathfrak{P}^{(2)}]^{}_{s}$ then,  in virtue of Lemma~\ref{LDCHNEchNHd}, we have that  $\mathfrak{P}'^{(2)}$ is an echelonless second-order path of length at least one that is not head-constant. Let $j\in\bb{\mathfrak{P}'^{(2)}}$ be the maximum index for which $\mathfrak{P}'^{(2),0,j}$ is a head-constant echelonless second-order path.

Thus, the value of the second-order Curry-Howard mapping at $\mathfrak{P}'^{(2)}$ is given by
$$
\mathrm{CH}^{(2)}_{s}\left(
\mathfrak{P}'^{(2)}
\right)
=
\mathrm{CH}^{(2)}_{s}\left(
\mathfrak{P}'^{(2),j+1,\bb{\mathfrak{P}'^{(2)}}-1}
\right)
\circ^{1\mathbf{T}_{\Sigma^{\boldsymbol{\mathcal{A}}^{(2)}}}(X)}_{s}
\mathrm{CH}^{(2)}_{s}\left(
\mathfrak{P}'^{(2),0,j}
\right).
$$

Now, we consider the different possibilities for $\mathfrak{Q}^{(2)}$. Note that, since $(\mathfrak{Q}^{(2)}, t)\prec_{\mathbf{Pth}_{\boldsymbol{\mathcal{A}}^{(2)}}} (\mathfrak{P}^{(2)},s)$, then following Definition~\ref{DDOrd}, either~(3.1.1) $\mathfrak{Q}^{(2)}=\mathfrak{P}^{(2),0,i}$ or~(3.1.2) $\mathfrak{Q}^{(2)}=\mathfrak{P}^{(2),i+1,\bb{\mathfrak{P}^{(2)}}-1}$. In any case, note that $t=s$.

If~(3.1.1), i.e., $\mathfrak{Q}^{(2)}=\mathfrak{P}^{(2),0,i}$ then, following Definition~\ref{DDOrd}, $(\mathfrak{P}'^{(2),0,j},s)\prec_{\mathbf{Pth}_{\boldsymbol{\mathcal{A}}^{(2)}}}(\mathfrak{P}'^{(2)},s)$. Moreover, since $\mathfrak{P}'^{(2)}$ is a second-order path in $[\mathfrak{P}^{(2)}]^{}_{s}$ we conclude that $\mathrm{CH}^{(2)}_{s}(\mathfrak{Q}^{(2)})=\mathrm{CH}^{(2)}_{s}(\mathfrak{P}'^{(2),0,j})$.

On the other hand, if~(3.1.2), i.e., $\mathfrak{Q}^{(2)}=\mathfrak{P}^{(2),i+1,\bb{\mathfrak{P}^{(2)}}-1}$ then, following Definition~\ref{DDOrd}, $(\mathfrak{P}'^{(2),j+1,\bb{\mathfrak{P}'^{(2)}}-1},s)\prec_{\mathbf{Pth}_{\boldsymbol{\mathcal{A}}^{(2)}}}(\mathfrak{P}'^{(2)},s)$. Moreover, since $\mathfrak{P}'^{(2)}$ is a second-order path in $[\mathfrak{P}^{(2)}]^{}_{s}$ we conclude that $\mathrm{CH}^{(2)}_{s}(\mathfrak{Q}^{(2)})=\mathrm{CH}^{(2)}_{s}(\mathfrak{P}'^{(2),j+1,\bb{\mathfrak{P}^{(2)}}-1})$.

The case of $\mathfrak{P}^{(2)}$ being not head-constant follows. 

We move to case~(3.2), where $\mathfrak{P}^{(2)}$ is a head-constant echelonless second-order path that is not coherent. Let $i\in\bb{\mathfrak{P}^{(2)}}$ be the maximum index for which $\mathfrak{P}^{(2),0,i}$ is a coherent head-constant echelonless second-order path. Following Definition~\ref{DDCH}, the value of the second-order Curry-Howard mapping at $\mathfrak{P}^{(2)}$ is given by
$$
\mathrm{CH}^{(2)}_{s}\left(
\mathfrak{P}^{(2)}
\right)
=
\mathrm{CH}^{(2)}_{s}\left(
\mathfrak{P}^{(2),i+1,\bb{\mathfrak{P}^{(2)}}-1}
\right)
\circ^{1\mathbf{T}_{\Sigma^{\boldsymbol{\mathcal{A}}^{(2)}}}(X)}_{s}
\mathrm{CH}^{(2)}_{s}\left(
\mathfrak{P}^{(2),0,i}
\right).
$$

Since $\mathfrak{P}^{(2)}$ is a head-constant echelonless second-order path of length at least one that is not coherent then, in virtue of Lemma~\ref{LDCHNEchHdNC}, we have that $\mathrm{CH}^{(2)}_{s}(\mathfrak{P}^{(2)})\in 
[\mathrm{T}_{\Sigma^{\boldsymbol{\mathcal{A}}^{(2)}}}(X)]^{\mathsf{HdC}\And\neg\mathsf{C}}_{s}
$.  Thus, if $\mathfrak{P}'^{(2)}$ is any second-order path in $[\mathfrak{P}^{(2)}]^{}_{s}$ then,  in virtue of Lemma~\ref{LDCHNEchHdNC}, we have that  $\mathfrak{P}'^{(2)}$ is a head-constant echelonless second-order path that is not coherent. Let $j\in\bb{\mathfrak{P}'^{(2)}}$ be the maximum index for which $\mathfrak{P}'^{(2),0,j}$ is a coherent head-constant echelonless second-order path.

Thus, the value of the second-order Curry-Howard mapping at $\mathfrak{P}'^{(2)}$ is given by
$$
\mathrm{CH}^{(2)}_{s}\left(
\mathfrak{P}'^{(2)}
\right)
=
\mathrm{CH}^{(2)}_{s}\left(
\mathfrak{P}'^{(2),j+1,\bb{\mathfrak{P}'^{(2)}}-1}
\right)
\circ^{1\mathbf{T}_{\Sigma^{\boldsymbol{\mathcal{A}}^{(2)}}}(X)}_{s}
\mathrm{CH}^{(2)}_{s}\left(
\mathfrak{P}'^{(2),0,j}
\right).
$$

Now, we consider the different possibilities for $\mathfrak{Q}^{(2)}$. Note that, since $(\mathfrak{Q}^{(2)}, t)\prec_{\mathbf{Pth}_{\boldsymbol{\mathcal{A}}^{(2)}}} (\mathfrak{P}^{(2)},s)$, then following Definition~\ref{DDOrd}, either~(3.2.1) $\mathfrak{Q}^{(2)}=\mathfrak{P}^{(2),0,i}$ or~(3.2.2) $\mathfrak{Q}^{(2)}=\mathfrak{P}^{(2),i+1,\bb{\mathfrak{P}^{(2)}}-1}$. In any case, note that $t=s$.

If~(3.2.1), i.e., if $\mathfrak{Q}^{(2)}=\mathfrak{P}^{(2),0,i}$ then, following Definition~\ref{DDOrd}, $(\mathfrak{P}'^{(2),0,j},s)\prec_{\mathbf{Pth}_{\boldsymbol{\mathcal{A}}^{(2)}}}(\mathfrak{P}'^{(2)},s)$. Moreover, since $\mathfrak{P}'^{(2)}$ is a second-order path in $[\mathfrak{P}^{(2)}]^{}_{s}$ we conclude that $\mathrm{CH}^{(2)}_{s}(\mathfrak{Q}^{(2)})=\mathrm{CH}^{(2)}_{s}(\mathfrak{P}'^{(2),0,j})$.

On the other hand, if~(3.2.2), i.e., if $\mathfrak{Q}^{(2)}=\mathfrak{P}^{(2),i+1,\bb{\mathfrak{P}^{(2)}}-1}$ then, following Definition~\ref{DDOrd}, $(\mathfrak{P}'^{(2),j+1,\bb{\mathfrak{P}'^{(2)}}-1},s)\prec_{\mathbf{Pth}_{\boldsymbol{\mathcal{A}}^{(2)}}}(\mathfrak{P}'^{(2)},s)$. Moreover, since $\mathfrak{P}'^{(2)}$ is a second-order path in $[\mathfrak{P}^{(2)}]^{}_{s}$ we conclude that $\mathrm{CH}^{(2)}_{s}(\mathfrak{Q}^{(2)})=\mathrm{CH}^{(2)}_{s}(\mathfrak{P}'^{(2),j+1,\bb{\mathfrak{P}^{(2)}}-1})$.

The case of $\mathfrak{P}^{(2)}$ being not coherent follows.

Finally, consider the case~(3.3), in which $\mathfrak{P}^{(2)}$ is a coherent head-constant echelonless second-order path.  Then, following Definition~\ref{DDHeadCt}, there exists a unique word $\mathbf{s}\in S^{\star}-\{\lambda\}$ and a unique operation symbol $\tau\in\Sigma^{\boldsymbol{\mathcal{A}}}_{\mathbf{s},s}$ associated to $\mathfrak{P}^{(2)}$. Let $(\mathfrak{P}^{(2)}_{j})_{j\in\bb{\mathbf{s}}}$ be the family of second-order paths we can extract from $\mathfrak{P}^{(2)}$ in virtue of Lemma~\ref{LDPthExtract}.  Then, following Definition~\ref{DDCH}, the value of the second-order Curry-Howard mapping at $\mathfrak{P}^{(2)}$ is given by
$$
\mathrm{CH}^{(2)}_{s}
\left(
\mathfrak{P}^{(2)}
\right)
=
\tau^{\mathbf{T}_{\Sigma^{\boldsymbol{\mathcal{A}}^{(2)}}}(X)}
\left(\left(\mathrm{CH}^{(2)}_{s_{j}}\left(
\mathfrak{P}^{(2)}_{j}
\right)\right)_{j\in\bb{\mathbf{s}}}\right).
$$

Since $\mathfrak{P}^{(2)}$ is a coherent echelonless second-order path associated to the operation symbol $\tau$ then, following  Lemma~\ref{LDCHNEchHdC}, we have that $\mathrm{CH}^{(2)}_{s}(\mathfrak{P}^{(2)})\in \mathcal{T}(\tau,\mathrm{T}_{\Sigma^{\boldsymbol{\mathcal{A}}}}(X))_{1}$, which is a subset of  $
[\mathrm{T}_{\Sigma^{\boldsymbol{\mathcal{A}}^{(2)}}}(X)]^{\mathsf{HdC}\And\mathsf{C}}_{s}
$.  Thus, if $\mathfrak{P}'^{(2)}$ is any second-order path in $[\mathfrak{P}^{(2)}]^{}_{s}$ then,  in virtue of Lemma~\ref{LDCHNEchHdC}, we have that  $\mathfrak{P}'^{(2)}$ is a coherent echelonless second-order path associated to the operation symbol $\tau$.

Let $(\mathfrak{P}'^{(2)}_{j})_{j\in\bb{\mathbf{s}}}$ be the family of second-order paths we can extract from $\mathfrak{P}'^{(2)}$ in virtue of Lemma~\ref{LDPthExtract}. Thus, the value of the second-order Curry-Howard mapping at $\mathfrak{P}'^{(2)}$ is given by
$$
\mathrm{CH}^{(2)}_{s}\left(
\mathfrak{P}'^{(2)}
\right)
=
\tau^{\mathbf{T}_{\Sigma^{\boldsymbol{\mathcal{A}}^{(2)}}}(X)}
\left(\left(\mathrm{CH}^{(2)}_{s_{j}}\left(
\mathfrak{P}'^{(2)}_{j}
\right)\right)_{j\in\bb{\mathbf{s}}}
\right).
$$

Now, we consider the different possibilities for $\mathfrak{Q}^{(2)}$. Note that, since $(\mathfrak{Q}^{(2)}, t)\prec_{\mathbf{Pth}_{\boldsymbol{\mathcal{A}}^{(2)}}} (\mathfrak{P}^{(2)},s)$ then following Definition~\ref{DDOrd}, there exist and index $j\in\bb{\mathbf{s}}$ for which $\mathfrak{Q}^{(2)}=\mathfrak{P}^{(2)}_{j}$. In this case $t=s_{j}$.

Then, following Definition~\ref{DDOrd},  $(\mathfrak{P}'^{(2)}_{j},s_{j})\prec_{\mathbf{Pth}_{\boldsymbol{\mathcal{A}}^{(2)}}}(\mathfrak{P}'^{(2)},s)$. Moreover, since $\mathfrak{P}'^{(2)}$ is a second-order path in $[\mathfrak{P}^{(2)}]^{}_{s}$ we conclude that $\mathrm{CH}^{(2)}_{s_{j}}(\mathfrak{Q}^{(2)})=\mathrm{CH}^{(2)}_{s_{j}}(\mathfrak{P}'^{(2)}_{j})$.

This completes the Case~(3).

This concludes the base case.

\textsf{Inductive step of the induction.}

Assume the statement holds for sequences of length up to $m$, i.e., for every pair of sorts $t,s\in S$, if $\mathfrak{Q}^{(2)}$ is a second-order path in $\mathrm{Pth}_{\boldsymbol{\mathcal{A}}^{(2)},t}$ and $\mathfrak{P}'^{(2)},\mathfrak{P}^{(2)}$ are second-order paths in $\mathrm{Pth}_{\boldsymbol{\mathcal{A}}^{(2)},s}$ such that  there exists a word $\mathbf{w}\in S^{\star}$ of length $\bb{\mathbf{w}}=m+1$ and a family of second-order paths $(\mathfrak{R}^{(2)}_{k})_{k\in\bb{\mathbf{w}}}$ in $\mathrm{Pth}_{\boldsymbol{\mathcal{A}}^{(2)},\mathbf{w}}$  such that $w_{0}=t$, $\mathfrak{R}^{(2)}_{0}=\mathfrak{Q}^{(2)}$, $w_{m}=s$, $\mathfrak{R}^{(2)}_{m}=\mathfrak{P}^{(2)}$ and 
for every $k\in m$, $(\mathfrak{R}^{(2)}_{k}, w_{k})\prec_{\mathbf{Pth}_{\boldsymbol{\mathcal{A}}^{(2)}}} (\mathfrak{R}^{(2)}_{k+1}, w_{k+1})$ and  $\mathfrak{P}'^{(2)}$ is a second-order path in $[\mathfrak{P}^{(2)}]^{}_{s}$, then there exists a second-order path $\mathfrak{Q}'^{(2)}$ in $[\mathfrak{Q}^{(2)}]^{}_{t}$ such that 
$$
\left(
\mathfrak{Q}'^{(2)},t
\right)\leq_{\mathbf{Pth}_{\boldsymbol{\mathcal{A}}^{(2)}}}
\left(
\mathfrak{P}'^{(2)},s
\right).
$$

Now we prove it for sequences of length $m+1$. 

Let $t,s$ be sorts in $S$, if $\mathfrak{Q}^{(2)}$ is a second-order path in $\mathrm{Pth}_{\boldsymbol{\mathcal{A}}^{(2)},t}$ and $\mathfrak{P}'^{(2)},\mathfrak{P}^{(2)}$ are second-order paths in $\mathrm{Pth}_{\boldsymbol{\mathcal{A}}^{(2)},s}$ such that  there exists a word $\mathbf{w}\in S^{\star}$ of length $\bb{\mathbf{w}}=m+2$ and a family of second-order paths $(\mathfrak{R}^{(2)}_{k})_{k\in\bb{\mathbf{w}}}$ in $\mathrm{Pth}_{\boldsymbol{\mathcal{A}}^{(2)},\mathbf{w}}$  such that $w_{0}=t$, $\mathfrak{R}^{(2)}_{0}=\mathfrak{Q}^{(2)}$, $w_{m+1}=s$, $\mathfrak{R}^{(2)}_{m+1}=\mathfrak{P}^{(2)}$ and 
for every $k\in m+1$, $(\mathfrak{R}^{(2)}_{k}, w_{k})\prec_{\mathbf{Pth}_{\boldsymbol{\mathcal{A}}^{(2)}}} (\mathfrak{R}^{(2)}_{k+1}, w_{k+1})$ and let  $\mathfrak{P}'^{(2)}$ be a second-order path in $[\mathfrak{P}^{(2)}]^{}_{s}$.

Consider the word $\mathbf{w}^{1,m+1}$ of length $\bb{\mathbf{w}^{1,m+1}}=m+1$ and the family of second-order paths $(\mathfrak{R}^{(2)}_{k})_{k\in\bb{\mathbf{w}^{1,m+1}}}$ in $\mathrm{Pth}_{\boldsymbol{\mathcal{A}}^{(2)},\mathbf{w}^{1,m+1}}$. This is a sequence of length $m$ instantiating that 
$
(\mathfrak{R}^{(2)}_{1},w_{1})
\leq_{\mathbf{Pth}_{\boldsymbol{\mathcal{A}}^{(2)}}}
(\mathfrak{P}^{(2)}, s).
$ Since $\mathfrak{P}'^{(2)}$ is a second-order path in $[\mathfrak{P}^{(2)}]^{}_{s}$, by the inductive hypothesis, there exist a second-order path $\mathfrak{R}'^{(2)}_{1}$ in $[\mathfrak{R}^{(2)}_{1}]^{}_{w_{1}}$ satisfying that 
$
(\mathfrak{R}'^{(2)}_{1},w_{1})
\leq_{\mathbf{Pth}_{\boldsymbol{\mathcal{A}}^{(2)}}}
(\mathfrak{P}'^{(2)},s).
$

Now, consider the sequence of second-order paths $(\mathfrak{Q}^{(2)},\mathfrak{R}^{(2)}_{1})$. It is a one-step sequence of second-order paths in $\mathrm{Pth}_{\boldsymbol{\mathcal{A}}^{(2)},\mathbf{w}^{0,1}}$ satisfying that $(\mathfrak{Q}^{(2)},t)\prec_{\mathbf{Pth}_{\boldsymbol{\mathcal{A}}^{(2)}}} (\mathfrak{R}^{(2)}_{1}, w_{1})$. Since $\mathfrak{R}'^{(2)}_{1}$ in $[\mathfrak{R}^{(2)}_{1}]^{}_{w_{1}}$, by the base case, we can find a second-order path $\mathfrak{Q}'^{(2)}$ in $[\mathfrak{Q}^{(2)}]^{}_{s}$ satisfying that 
$
(\mathfrak{Q}'^{(2)},t)
\leq_{\mathbf{Pth}_{\boldsymbol{\mathcal{A}}^{(2)}}}
(\mathfrak{R}'^{(2)}_{1}, w_{1}).
$	

Hence, $\mathfrak{Q}'^{(2)}$ is a second-order path in $[\mathfrak{Q}^{(2)}]^{}_{s}$ satisfying that 
$$
\left(
\mathfrak{Q}'^{(2)},t
\right)
\leq_{\mathbf{Pth}_{\boldsymbol{\mathcal{A}}^{(2)}}}
\left(
\mathfrak{P}'^{(2)}, s
\right).
$$

This concludes the proof.
\end{proof}

\begin{corollary}\label{CDCHOrdII} Let $t,s$ be sorts in $S$, $[\mathfrak{Q}^{(2)}]^{}_{t}$ a second-order path class in $[\mathrm{Pth}_{\boldsymbol{\mathcal{A}}^{(2)}}]_{t}$ and $[\mathfrak{P}^{(2)}]^{}_{s}$ a second-order path class in $[\mathrm{Pth}_{\boldsymbol{\mathcal{A}}^{(2)}}]_{s}$. Then the following statements are equivalent
\begin{enumerate}
\item[(i)] $([\mathfrak{Q}^{(2)}]^{}_{t},t)
\leq_{[\mathbf{Pth}_{\boldsymbol{\mathcal{A}}^{(2)}}]}
([\mathfrak{P}^{(2)}]^{}_{s},s);
$
\item[(ii)] There exists $
\mathfrak{Q}'^{(2)}\in [\mathfrak{Q}^{(2)}]^{}_{t}$ such that $
(\mathfrak{Q}'^{(2)},t)
\leq_{\mathbf{Pth}_{\boldsymbol{\mathcal{A}}^{(2)}}}
(\mathfrak{P}^{(2)},s)
.
$
\end{enumerate}
\end{corollary}
\begin{proof}
Assume that $([\mathfrak{Q}^{(2)}]^{}_{t},t)
\leq_{[\mathbf{Pth}_{\boldsymbol{\mathcal{A}}^{(2)}}]}
([\mathfrak{P}^{(2)}]^{}_{s},s)
$ then there exists $\mathfrak{Q}'^{(2)}\in[\mathfrak{Q}^{(2)}]^{}_{t}$ and $\mathfrak{P}'^{(2)}\in [\mathfrak{P}^{(2)}]^{}_{s}$ such that $(\mathfrak{Q}'^{(2)},t) \leq_{\mathbf{Pth}_{\boldsymbol{\mathcal{A}}^{(2)}}} (\mathfrak{P}'^{(2)},s)$. Since $\mathfrak{P}^{(2)}\in[\mathfrak{P}'^{(2)}]^{}_{s}$, by Lemma~\ref{LDCHOrdII}, we can find $\mathfrak{Q}''^{(2)}\in[\mathfrak{Q}'^{(2)}]^{}_{t}$ such that $(\mathfrak{Q}''^{(2)},t) \leq_{\mathbf{Pth}_{\boldsymbol{\mathcal{A}}^{(2)}}} (\mathfrak{P}^{(2)},s)$. Note that $\mathfrak{Q}''^{(2)}\in[\mathfrak{Q}^{(2)}]^{}_{t}$.

The other implication follows by definition of the relation $\leq_{[\mathbf{Pth}_{\boldsymbol{\mathcal{A}}^{(2)}}]}$.

This concludes the proof.
\end{proof}

Next lemma shows that if we are given two second-order paths in the same class under the Kernel of the second-order Curry Howard mapping then these second-order paths can only be related with respect to the partial order $\leq_{\mathbf{Pth}_{\boldsymbol{\mathcal{A}}^{(2)}}}$ if they are equal.

\begin{lemma}\label{LDCHOrdIII} Let $s$ be a sort in $S$ and $\mathfrak{P}'^{(2)},\mathfrak{P}^{(2)}$ second-order paths in $\mathrm{Pth}_{\boldsymbol{\mathcal{A}}^{(2)},s}$. If  $(\mathfrak{P}'^{(2)}, \mathfrak{P}^{(2)})\in\mathrm{Ker}(\mathrm{CH}^{(2)})_{s}$ and $(\mathfrak{P}'^{(2)},s)\leq_{\mathbf{Pth}_{\boldsymbol{\mathcal{A}}^{(2)}}}(\mathfrak{P}^{(2)},s)$ then $\mathfrak{P}'^{(2)}=\mathfrak{P}^{(2)}$.
\end{lemma}
\begin{proof}
Assume towards a contradiction that $\mathfrak{P}'^{(2)}\neq\mathfrak{P}^{(2)}$. Since $(\mathfrak{P}'^{(2)},s)\leq_{\mathbf{Pth}_{\boldsymbol{\mathcal{A}}^{(2)}}}(\mathfrak{P}^{(2)},s)$, we conclude that $(\mathfrak{P}'^{(2)},s)$ is strictly smaller than $(\mathfrak{P}^{(2)},s)$ with respect to the partial order $\leq_{\mathbf{Pth}_{\boldsymbol{\mathcal{A}}^{(2)}}}$. Following Proposition~\ref{PDCHMono} this entails that $(\mathrm{CH}^{(2)}_{s}(\mathfrak{P}'^{(2)}),s)$ is strictly smaller than $(\mathrm{CH}^{(2)}_{s}(\mathfrak{P}^{(2)}),s)$  with respect to the partial order $\leq_{\mathbf{T}_{\Sigma^{\boldsymbol{\mathcal{A}}^{(2)}}}(X)}$, contradicting the fact that, by hypothesis, $(\mathfrak{P}'^{(2)}, \mathfrak{P}^{(2)})\in\mathrm{Ker}(\mathrm{CH}^{(2)})_{s}$.
Therefore, $\mathfrak{P}'^{(2)}=\mathfrak{P}^{(2)}$.
\end{proof}

As a corollary of the previous lemma we have that no pair can be found with respect to the partial order $\leq_{\mathbf{Pth}_{\boldsymbol{\mathcal{A}}^{(2)}}}$ between two equivalent second-order paths.

\begin{corollary}\label{CDCHOrdIII} Let $s,t$ be sorts in $S$, $\mathfrak{Q}^{(2)}$ a second-order path in $\mathrm{Pth}_{\boldsymbol{\mathcal{A}}^{(2)},t}$ and $\mathfrak{P}'^{(2)},\mathfrak{P}^{(2)}$  two second-order paths in $\mathrm{Pth}_{\boldsymbol{\mathcal{A}}^{(2)},s}$. If 
$$
\left(
\mathfrak{P}'^{(2)},s
\right)
\leq_{\mathbf{Pth}_{\boldsymbol{\mathcal{A}}^{(2)}}}
\left(
\mathfrak{Q}^{(2)},t
\right)
\leq_{\mathbf{Pth}_{\boldsymbol{\mathcal{A}}^{(2)}}}
\left(
\mathfrak{P}^{(2)},s
\right)
$$
and  $(\mathfrak{P}'^{(2)},\mathfrak{P}^{(2)})\in\mathrm{Ker}(\mathrm{CH}^{(2)})_{s}$, then 
$t=s$ and $\mathfrak{Q}^{(2)}=\mathfrak{P}^{(2)}$.
\end{corollary}
\begin{proof}
According to Lemma~\ref{LDCHOrdIII}, $\mathfrak{P}'^{(2)}=\mathfrak{P}^{(2)}$. The statement follows from the antisymmetry of the preorder $\leq_{\mathbf{Pth}_{\boldsymbol{\mathcal{A}}^{(2)}}}$ proven in Proposition~\ref{PDOrdArt}.
\end{proof}

We are now in position to prove that $\leq_{[\mathbf{Pth}_{\boldsymbol{\mathcal{A}}^{(2)}}]}$ is an Artinian partial order defined on the set $\coprod[\mathrm{Pth}_{\boldsymbol{\mathcal{A}}^{(2)}}]$.

\begin{restatable}{proposition}{PDCHOrdArt}
\label{PDCHOrdArt}  $(\coprod[\mathrm{Pth}_{\boldsymbol{\mathcal{A}}^{(2)}}], \leq_{[\mathbf{Pth}_{\boldsymbol{\mathcal{A}}^{(2)}}]})$ is a partially ordered set. Moreover, in this partially ordered set there is not any strictly decreasing $\omega_{0}$-chain, i.e.,  $(\coprod[\mathrm{Pth}_{\boldsymbol{\mathcal{A}}^{(2)}}], \leq_{[\mathbf{Pth}_{\boldsymbol{\mathcal{A}}^{(2)}}]})$ is an Artinian partially ordered set.
\end{restatable}

\begin{proof}
That $\leq_{[\mathbf{Pth}_{\boldsymbol{\mathcal{A}}^{(2)}}]}$ is reflexive follows from the fact that $\leq_{\mathbf{Pth}_{\boldsymbol{\mathcal{A}}^{(2)}}}$ is reflexive. Indeed, let $s$ be a sort in $S$ and let $[\mathfrak{P}^{(2)}]^{}_{s}$ be a second-order path class in $\mathrm{Pth}_{\boldsymbol{\mathcal{A}}^{(2)},s}$, then since $\mathfrak{P}^{(2)}\in [\mathfrak{P}^{(2)}]^{}_{s}$ and $(\mathfrak{P}^{(2)},s)\leq_{\mathbf{Pth}_{\boldsymbol{\mathcal{A}}^{(2)}}} (\mathfrak{P}^{(2)},s)$, then we conclude that 
$$
\left(
\left[\mathfrak{P}^{(2)}
\right]^{}_{s},s
\right)
\leq_{[\mathbf{Pth}_{\boldsymbol{\mathcal{A}}^{(2)}}]} 
\left(
\left[
\mathfrak{P}^{(2)}
\right]^{}_{s},s
\right).$$

We now prove that $\leq_{[\mathbf{Pth}_{\boldsymbol{\mathcal{A}}^{(2)}}]}$ is antisymmetric. Let $t,s$ be sorts in $S$ and let $[\mathfrak{Q}^{(2)}]^{}_{t}$ be a second-order path class in $[\mathrm{Pth}_{\boldsymbol{\mathcal{A}}^{(2)}}]_{t}$ and let $[\mathfrak{P}^{(2)}]^{}_{s}$ be a second-order path class in $[\mathrm{Pth}_{\boldsymbol{\mathcal{A}}^{(2)}}]_{s}$ satisfying that
\begin{itemize}
\item[(i)] $\left(\left[\mathfrak{Q}^{(2)}
\right]^{}_{t},t
\right)
\leq_{[\mathbf{Pth}_{\boldsymbol{\mathcal{A}}^{(2)}}]}\left(
\left[\mathfrak{P}^{(2)}
\right]^{}_{s},s\right);$
\item[(ii)] $\left(
\left[
\mathfrak{P}^{(2)}
\right]^{}_{s},s
\right)
\leq_{[\mathbf{Pth}_{\boldsymbol{\mathcal{A}}^{(2)}}]}\left(\left[\mathfrak{Q}^{(2)}
\right]^{}_{t},t\right).$
\end{itemize}

We want to prove that $t=s$ and $[\mathfrak{Q}^{(2)}]^{}_{t}=[\mathfrak{P}^{(2)}]^{}_{s}$.

On one hand, since the inequality (i) holds, we can find second-order paths $\mathfrak{Q}'^{(2)}\in[\mathfrak{Q}^{(2)}]^{}_{t}$ and $\mathfrak{P}'^{(2)}\in[\mathfrak{P}^{(2)}]^{}_{s}$ satisfying that 
$
(\mathfrak{Q}'^{(2)},t)\leq_{\mathbf{Pth}_{\boldsymbol{\mathcal{A}}^{(2)}}}
(\mathfrak{P}'^{(2)},s).
$ On the other hand, since the inequality (ii) holds, we can find second-order paths  $\mathfrak{P}''^{(2)}\in[\mathfrak{P}^{(2)}]^{}_{s}$ and $\mathfrak{Q}''^{(2)}\in[\mathfrak{Q}^{(2)}]^{}_{t}$ satisfying that 
$
(\mathfrak{P}''^{(2)},s)\leq_{\mathbf{Pth}_{\boldsymbol{\mathcal{A}}^{(2)}}}
(\mathfrak{Q}''^{(2)},t).
$

Since $
(\mathfrak{Q}'^{(2)},t)\leq_{\mathbf{Pth}_{\boldsymbol{\mathcal{A}}^{(2)}}}
(\mathfrak{P}'^{(2)},s)
$ and $\mathfrak{Q}''^{(2)}$ is a second-order path in $[\mathfrak{Q}'^{(2)}]^{}_{t}$, in virtue of Lemma~\ref{LDCHOrdI}, we can find a second-order path $\mathfrak{P}'''^{(2)}\in[\mathfrak{P}'^{(2)}]^{}_{s}$ for which
$
(\mathfrak{Q}''^{(2)},t)\leq_{\mathbf{Pth}_{\boldsymbol{\mathcal{A}}^{(2)}}}
(\mathfrak{P}'''^{(2)},s).
$

Hence, we have found the chain
$$
\left(
\mathfrak{P}''^{(2)},s
\right)
\leq_{\mathbf{Pth}_{\boldsymbol{\mathcal{A}}^{(2)}}}
\left(
\mathfrak{Q}''^{(2)},t
\right)
\leq_{\mathbf{Pth}_{\boldsymbol{\mathcal{A}}^{(2)}}}
\left(\mathfrak{P}'''^{(2)},s
\right).
$$

Moreover, $(\mathfrak{P}''^{(2)},\mathfrak{P}'''^{(2)})\in\mathrm{Ker}(\mathrm{CH}^{(2)})_{s}$. Hence, in virtue of Corollary~\ref{CDCHOrdIII}, we conclude that $t=s$ and $\mathfrak{Q}^{(2)}=\mathfrak{P}^{(2)}$, which implies that $[\mathfrak{Q}^{(2)}]^{}_{t}=[\mathfrak{P}^{(2)}]^{}_{s}$.

We now prove that $\leq_{[\mathbf{Pth}_{\boldsymbol{\mathcal{A}}^{(2)}}]}$ is transitive. Let $u,t,s$ be sorts in $S$ and let $[\mathfrak{R}^{(2)}]_{u}$ be a second-order path class $[\mathrm{Pth}_{\boldsymbol{\mathcal{A}}^{(2)}}]_{u}$,  $[\mathfrak{Q}^{(2)}]^{}_{t}$ be a second-order path class in $[\mathrm{Pth}_{\boldsymbol{\mathcal{A}}^{(2)}}]_{t}$ and let $[\mathfrak{P}^{(2)}]^{}_{s}$ be a second-order path class in $[\mathrm{Pth}_{\boldsymbol{\mathcal{A}}^{(2)}}]_{s}$ satisfying that
\begin{itemize}
\item[(i)] $\left(\left[
\mathfrak{R}^{(2)}
\right]_{u},u
\right)
\leq_{[\mathbf{Pth}_{\boldsymbol{\mathcal{A}}^{(2)}}]}\left(\left[
\mathfrak{Q}^{(2)}
\right]^{}_{t},t\right);$
\item[(ii)] $\left(
\left[
\mathfrak{Q}^{(2)}
\right]^{}_{t},t
\right)
\leq_{[\mathbf{Pth}_{\boldsymbol{\mathcal{A}}^{(2)}}]}
\left(\left[
\mathfrak{P}^{(2)}
\right]^{}_{s},s\right).$
\end{itemize}

We want to prove that $([\mathfrak{R}^{(2)}]_{u},u)\leq_{[\mathbf{Pth}_{\boldsymbol{\mathcal{A}}^{(2)}}]}([\mathfrak{P}^{(2)}]^{}_{s},s)$.

On one hand, since inequality (i) holds, we can find second-order paths $\mathfrak{R}'^{(2)}\in[\mathfrak{R}^{(2)}]_{u}$ and $\mathfrak{Q}'^{(2)}\in[\mathfrak{Q}^{(2)}]^{}_{t}$ satisfying that 
$
(\mathfrak{R}'^{(2)},u)\leq_{\mathbf{Pth}_{\boldsymbol{\mathcal{A}}^{(2)}}}
(\mathfrak{Q}'^{(2)},t).
$ On the other hand, since inequality (ii) holds, we can find second-order paths  $\mathfrak{Q}''^{(2)}\in[\mathfrak{Q}^{(2)}]^{}_{t}$ and $\mathfrak{P}'^{(2)}\in[\mathfrak{P}^{(2)}]^{}_{s}$ satisfying that 
$
(\mathfrak{Q}''^{(2)},t)\leq_{\mathbf{Pth}_{\boldsymbol{\mathcal{A}}^{(2)}}}
(\mathfrak{P}'^{(2)},s).
$

Since $
(\mathfrak{Q}''^{(2)},t)\leq_{\mathbf{Pth}_{\boldsymbol{\mathcal{A}}^{(2)}}}
(\mathfrak{P}'^{(2)},s)
$ and $\mathfrak{Q}'^{(2)}$ is a second-order path in $[\mathfrak{Q}''^{(2)}]^{}_{t}$, in virtue of Lemma~\ref{LDCHOrdI}, we can find a path $\mathfrak{P}''^{(2)}\in[\mathfrak{P}'^{(2)}]^{}_{s}$ for which
$
(\mathfrak{Q}''^{(2)},t)\leq_{\mathbf{Pth}_{\boldsymbol{\mathcal{A}}^{(2)}}}
(\mathfrak{P}''^{(2)},s).
$

Hence, we have found the chain
$$
\left(
\mathfrak{R}'^{(2)},u
\right)
\leq_{\mathbf{Pth}_{\boldsymbol{\mathcal{A}}^{(2)}}}
\left(
\mathfrak{Q}''^{(2)},t
\right)
\leq_{\mathbf{Pth}_{\boldsymbol{\mathcal{A}}^{(2)}}}
\left(
\mathfrak{P}''^{(2)},s
\right).
$$
 
By the transitivity of $\leq_{\mathbf{Pth}_{\boldsymbol{\mathcal{A}}^{(2)}}}$, we conclude that $(\mathfrak{R}'^{(2)},u)
\leq_{\mathbf{Pth}_{\boldsymbol{\mathcal{A}}^{(2)}}}
(\mathfrak{P}''^{(2)},s)$. Note that $\mathfrak{R}'^{(2)}\in[\mathfrak{R}^{(2)}]_{u}$ and $\mathfrak{P}''^{(2)}\in[\mathfrak{P}^{(2)}]^{}_{s}$. Hence, we have that 
$$
\left(\left[
\mathfrak{R}^{(2)}
\right]^{}_{u},u
\right)\leq_{[\mathbf{Pth}_{\boldsymbol{\mathcal{A}}^{(2)}}]}\left(
\left[
\mathfrak{P}^{(2)}
\right]^{}_{s},s
\right).$$

From this it follows that $\leq_{[\mathbf{Pth}_{\boldsymbol{\mathcal{A}}^{(2)}}]}$ is a partial order.

We next prove that no strictly decreasing $\omega_{0}$-chain can exist in the partially ordered set $(\coprod[\mathrm{Pth}_{\boldsymbol{\mathcal{A}}^{(2)}}],\leq_{[\mathbf{Pth}_{\boldsymbol{\mathcal{A}}^{(2)}}]})$. Let us suppose, towards a contradiction, that there exists at least one strictly decreasing $\omega_{0}$-chain $(([\mathfrak{R}^{(2)}_{k}]^{}_{s_{k}},s_{k}))_{k\in\omega_{0}}$ in $(\coprod[\mathrm{Pth}_{\boldsymbol{\mathcal{A}}^{(2)}}],\leq_{[\mathbf{Pth}_{\boldsymbol{\mathcal{A}}^{(2)}}]})$. 

Then, for every $k\in\omega_{0}$, we have that
\begin{itemize}
\item[(i)] $[\mathfrak{R}^{(2)}_{k}
]^{}_{s_{k}}$ is a second-order path class in  $[\mathrm{Pth}_{\boldsymbol{\mathcal{A}}^{(2)}}]_{s_{k}}$, and
\item[(ii)] $([
\mathfrak{R}^{(2)}_{k+1}
]^{}_{s_{k+1}}, s_{k+1}
)<_{[\mathbf{Pth}_{\boldsymbol{\mathcal{A}}^{(2)}}]} 
([
\mathfrak{R}^{(2)}_{k}
]^{}_{s_{k}}, s_{k}
).$
\end{itemize}

We next prove how this strictly decreasing $\omega_{0}$-chain in the partially ordered set $(\coprod[\mathrm{Pth}_{\boldsymbol{\mathcal{A}}^{(2)}}],\leq_{[\mathbf{Pth}_{\boldsymbol{\mathcal{A}}^{(2)}}]})$ leads to a strictly decreasing $\omega_{0}$-chain in the partially ordered set $(\coprod\mathrm{Pth}_{\boldsymbol{\mathcal{A}}^{(2)}},\leq_{\mathbf{Pth}_{\boldsymbol{\mathcal{A}}^{(2)}}})$. 

Let us define $\mathfrak{R}'^{(2)}_{0}$ to be equal to $\mathfrak{R}^{(2)}_{0}$. Since $([\mathfrak{R}^{(2)}_{1}]^{}_{s_{1}}, s_{1})<_{[\mathbf{Pth}_{\boldsymbol{\mathcal{A}}^{(2)}}]} ([\mathfrak{R}'^{(2)}_{0}]^{}_{s_{0}}, s_{0})$ then, in virtue of Corollary~\ref{CDCHOrdII}, we can find $\mathfrak{R}'^{(2)}_{1}\in[\mathfrak{R}^{(2)}_{1}]^{}_{s_{1}}$ satisfying that $(\mathfrak{R}'^{(2)}_{1},s_{1})$ $<_{\mathbf{Pth}_{\boldsymbol{\mathcal{A}}^{(2)}}}$-precedes $(\mathfrak{R}'^{(2)}_{0},s_{0})$. 

Following this procedure, assume that, for $n\in\mathbb{N}$ with $n\geq 2$, we have found a strictly-decreasing $n$-chain $((\mathfrak{R}'^{(2)}_{k},s_{k}))_{k\in n}$ in $(\coprod\mathrm{Pth}_{\boldsymbol{\mathcal{A}}^{(2)}},\leq_{\mathbf{Pth}_{\boldsymbol{\mathcal{A}}^{(2)}}})$ satisfying that, for every $k\in n$, $\mathfrak{R}'^{(2)}_{k}\in[\mathfrak{R}^{(2)}_{k}]^{}_{s_{k}}$. Now, since $([\mathfrak{R}^{(2)}_{n}]^{}_{s_{n}}, s_{n})<_{[\mathbf{Pth}_{\boldsymbol{\mathcal{A}}^{(2)}}]} ([\mathfrak{R}^{(2)}_{n-1}]^{}_{s_{n-1}}, s_{n-1})$ and $\mathfrak{R}'^{(2)}_{n-1}\in [\mathfrak{R}^{(2)}_{n-1}]^{}_{s_{n-1}}$ then, in virtue of Corollary~\ref{CDCHOrdII}, we can find $\mathfrak{R}'^{(2)}_{n}\in[\mathfrak{R}^{(2)}_{n}]^{}_{s_{n}}$ satisfying that $(\mathfrak{R}'^{(2)}_{n},s_{n})<_{\mathbf{Pth}_{\boldsymbol{\mathcal{A}}^{(2)}}}(\mathfrak{R}'^{(2)}_{n-1},s_{n-1})$.  Note that  $((\mathfrak{R}'^{(2)}_{k},s_{k}))_{k\in n+1}$ is  a strictly-decreasing $n+1$-chain in $(\coprod\mathrm{Pth}_{\boldsymbol{\mathcal{A}}^{(2)}},\leq_{\mathbf{Pth}_{\boldsymbol{\mathcal{A}}^{(2)}}})$ satisfying that, for every $k\in n+1$, $\mathfrak{R}'^{(2)}_{k}\in[\mathfrak{R}^{(2)}_{k}]^{}_{s_{k}}$.

This ultimately leads to a strictly-decreasing $\omega_{0}$-chain in $(\coprod\mathrm{Pth}_{\boldsymbol{\mathcal{A}}^{(2)}},\leq_{\mathbf{Pth}_{\boldsymbol{\mathcal{A}}^{(2)}}})$ contradicting the fact that, by Proposition~\ref{PDOrdArt}, $(\coprod\mathrm{Pth}_{\boldsymbol{\mathcal{A}}^{(2)}},\leq_{\mathbf{Pth}_{\boldsymbol{\mathcal{A}}^{(2)}}})$ is an Artinian partially ordered set.

It follows that $(\coprod[\mathrm{Pth}_{\boldsymbol{\mathcal{A}}^{(2)}}],\leq_{[\mathbf{Pth}_{\boldsymbol{\mathcal{A}}^{(2)}}]})$ is an Artinian partially ordered set.

This concludes the proof.
\end{proof}

\begin{restatable}{proposition}{PDCHOrd}
\label{PDCHOrd} The mapping $\coprod\mathrm{pr}^{\mathrm{Ker}(\mathrm{CH}^{(2)})}\colon\coprod\mathrm{Pth}_{\boldsymbol{\mathcal{A}}^{(2)}}\mor\coprod[\mathrm{Pth}_{\boldsymbol{\mathcal{A}}^{(2)}}]$ that, for every $s\in S$ maps a pair $(\mathfrak{P}^{(2)},s)$ in $\coprod\mathrm{Pth}_{\boldsymbol{\mathcal{A}}^{(2)}}$ to the pair $([\mathfrak{P}^{(2)}]^{}_{s},s)$ in $\coprod[\mathrm{Pth}_{\boldsymbol{\mathcal{A}}^{(2)}}]$ is order-preserving and inversely order-preserving.
\end{restatable}

Taking into account the definitions above, now the coproduct of the $(2,[1])$-identity second-order path mapping becomes an order embedding from the Artinian partial order $(\coprod[\mathrm{PT}_{\boldsymbol{\mathcal{A}}}], \leq_{[\mathbf{PT}_{\boldsymbol{\mathcal{A}}}]})$, introduced on Definition~\ref{DPTQOrd}, to the Artinian partial order $(\coprod[\mathrm{Pth}_{\boldsymbol{\mathcal{A}}^{(2)}}], \leq_{[\mathbf{Pth}_{\boldsymbol{\mathcal{A}}^{(2)}}]})$, introduced on Definition~\ref{DDCHOrd}.

\begin{restatable}{proposition}{PDCHOrdIp}
\label{PDCHOrdIp} 
The mapping $\coprod\mathrm{ip}^{([2],[1])\sharp}$ is an order-embedding
\[
\textstyle
\coprod\mathrm{ip}^{([2],[1])\sharp}\colon
\left(
\coprod
\left[
\mathrm{PT}_{\boldsymbol{\mathcal{A}}}
\right],
\leq_{[\mathbf{PT}_{\boldsymbol{\mathcal{A}}}]}
\right)
\mor 
\textstyle
\left(
\coprod[\mathrm{Pth}_{\boldsymbol{\mathcal{A}}^{(2)}}], \leq_{[\mathbf{Pth}_{\boldsymbol{\mathcal{A}}^{(2)}}]}
\right).
\]
\end{restatable}
\begin{proof}
Let $s,t$ be sorts in $S$ and $([Q]_{t},t), ([P]_{s},s)$  pairs in $\coprod
[\mathrm{PT}_{\boldsymbol{\mathcal{A}}}]$. We need to prove that the following statements are equivalent
\begin{enumerate}
\item[(i)] $([Q]_{t},t)
\leq_{[\mathbf{PT}_{\boldsymbol{\mathcal{A}}}]}
([P]_{s},s)$;
\item[(ii)] $(\mathrm{ip}^{([2],[1])\sharp}_{t}([Q]_{t}),t)
\leq_{[\mathbf{Pth}_{\boldsymbol{\mathcal{A}}^{(2)}}]}
(\mathrm{ip}^{([2],[1])\sharp}_{s}([P]_{s}),s).
$ 
\end{enumerate}

Note that the following chain of equivalences holds
\begin{flushleft}
$
\left([Q]_{t},t
\right)
\leq_{[\mathbf{PT}_{\boldsymbol{\mathcal{A}}}]}
\left(
[P]_{s},s
\right)$
\allowdisplaybreaks
\begin{align*}
\quad
&\Longleftrightarrow
\left(
\mathrm{ip}^{(2,[1])\sharp}_{t}
\left(
[Q]_{t}
\right),t
\right)
\leq_{\mathbf{Pth}_{\boldsymbol{\mathcal{A}}^{(2)}}}
\left(
\mathrm{ip}^{(2,[1])\sharp}_{s}
\left(
[P]_{s}
\right),s
\right)
\tag{1}
\\&\Longleftrightarrow
\left(
\left[
\mathrm{ip}^{(2,[1])\sharp}_{t}
\left(
[Q]_{t}
\right)\right]^{}_{t},t
\right)
\leq_{[\mathbf{Pth}_{\boldsymbol{\mathcal{A}}^{(2)}}]}
\left(
\left[
\mathrm{ip}^{(2,[1])\sharp}_{s}
\left(
[P]_{s}
\right)\right]^{}_{s},s\right)
\tag{2}
\\&\Longleftrightarrow
\left(
\mathrm{ip}^{([2],[1])\sharp}_{t}
\left(
[Q]_{t}
\right),t
\right)
\leq_{[\mathbf{Pth}_{\boldsymbol{\mathcal{A}}^{(2)}}]}
\left(
\mathrm{ip}^{([2],[1])\sharp}_{s}
\left(
[P]_{s}
\right),s\right).
\tag{3}
\end{align*}
\end{flushleft}

The first equivalence follows from the fact that, by Proposition~\ref{PDUOrdEmb}, the mapping $\coprod\mathrm{ip}^{(2,[1])\sharp}$ is an order embedding; the second equivalence follows from left to right from the fact that, by Proposition~\ref{PDCHOrd}, $\coprod\mathrm{pr}^{\mathrm{Ker}(\mathrm{CH}^{(2)})}$ is order-preserving. On the other hand the implication from right to left follows from the fact that, by Proposition~\ref{PDCHOrd}, $\coprod\mathrm{pr}^{\mathrm{Ker}(\mathrm{CH}^{(2)})}$ is inversely order-preserving and from the fact that, by Corollary~\ref{CDCHUId}, the class under the kernel of the second-order Curry-Howard mapping of a $(2,[1])$-identity second-order path only contains this exact $(2,[1])$-identity second-order path; finally, the last equivalence follows from the fact that, according to Definition~\ref{DDCHDU}, $\mathrm{ip}^{([2],[1])\sharp}=\mathrm{pr}^{\mathrm{Ker}(\mathrm{CH}^{(2)})}\circ\mathrm{ip}^{(2,[1])\sharp}$.

This concludes the proof.
\end{proof}

\begin{restatable}{proposition}{PDCHOrdMono}
\label{PDCHOrdMono} 
The mapping $\coprod\mathrm{CH}^{(2)\mathrm{m}}$ from
$\coprod[\mathrm{Pth}_{\boldsymbol{\mathcal{A}}^{(2)}}]$ to  $\coprod\mathrm{T}_{\Sigma^{\boldsymbol{\mathcal{A}}^{(2)}}}(X)$ determines an order-preserving mapping 
\[
\textstyle
\coprod\mathrm{CH}^{(2)\mathrm{m}}\colon
\left(
\coprod[\mathrm{Pth}_{\boldsymbol{\mathcal{A}}^{(2)}}], \leq_{[\mathbf{Pth}_{\boldsymbol{\mathcal{A}}^{(2)}}]}
\right)
\mor
\textstyle
\left(
\coprod\mathrm{T}_{\Sigma^{\boldsymbol{\mathcal{A}}^{(2)}}}(X), \leq_{\mathbf{T}_{\Sigma^{\boldsymbol{\mathcal{A}}^{(2)}}}(X)}\right),
\]
i.e., given pairs $([\mathfrak{Q}^{(2)}]^{}_{t},t),([\mathfrak{P}^{(2)}]^{}_{s},s)$ in $\coprod[\mathrm{Pth}_{\boldsymbol{\mathcal{A}}^{(2)}}]$ satisfying that  
$$
\left(\left[
\mathfrak{Q}^{(2)}
\right]^{}_{t},t
\right)
\leq_{[\mathbf{Pth}_{\boldsymbol{\mathcal{A}}^{(2)}}]} 
\left(\left[
\mathfrak{P}^{(2)}
\right]^{}_{s},s\right),$$

then $$\left(
\mathrm{CH}^{(2)\mathrm{m}}_{t}
\left(\left[
\mathfrak{Q}^{(2)}
\right]^{}_{t}
\right),t\right)
\leq_{\mathbf{T}_{\Sigma^{\boldsymbol{\mathcal{A}}^{(2)}}}(X)}
\left(
\mathrm{CH}^{(2)\mathrm{m}}_{s}
\left(\left[
\mathfrak{P}^{(2)}
\right]^{}_{s}
\right),s
\right),$$
that is, $\mathrm{CH}^{(2)}_{t}(\mathfrak{Q}^{(2)})$ is a subterm of type $t$ of the term $\mathrm{CH}^{(2)}_{s}(\mathfrak{P}^{(2)})$.
\end{restatable}
\begin{proof}
It follows from Proposition~\ref{PDCHMono}
\end{proof}

\chapter{
\texorpdfstring
{The congruence $\Upsilon^{[1]}$ on $\mathbf{Pth}_{\boldsymbol{\mathcal{A}}^{(2)}}$ and its relatives}
{The Upsilon congruence}
}\label{S2G}

In this chapter we begin by defining the relation $\Upsilon^{(1)}$ on $\mathrm{Pth}_{\boldsymbol{\mathcal{A}}^{(2)}}$ with the aim of recovering the categorial identities on the operations of $0$-composition, $0$-source and $0$-target. We first prove that pairs of second-order paths in $\Upsilon^{(1)}$ have the same length, the same $([1],2)$-source and $([1],2)$-target. As a consequence, we will show that $(2,[1])$-identity second-order paths in $\Upsilon^{(1)}$ are equal. We also consider the smallest $\Sigma^{\boldsymbol{\mathcal{A}}^{(2)}}$-congruence containing $\Upsilon^{(1)}$, denoted by $\Upsilon^{[1]}$ and we check that pairs of second-order paths in $\Upsilon^{[1]}$ have the same length, the same $([1],2)$-source and $([1],2)$-target. As a result, pairs of second-order paths in $\Upsilon^{[1]}$ also have the same $(0,2)$-source and $(0,2)$-target. As before, we will check that $(2,[1])$-identity second-order paths in $\Upsilon^{[1]}$ are equal. We will denote by $[\mathbf{Pth}_{\boldsymbol{\mathcal{A}}^{(2)}}]_{\Upsilon^{[1]}}$ to the many-sorted partial $\Sigma^{\boldsymbol{\mathcal{A}}^{(2)}}$-algebra on the quotient. We conclude this chapter by looking at how to order this quotient. We define the relation $\leq_{[\mathbf{Pth}_{\boldsymbol{\mathcal{A}}^{(2)}}]_{\Upsilon^{[1]}}}$ and prove that it is an Artinian partial preorder for which $\mathrm{pr}^{\Upsilon^{[1]}}$, the projection mapping to the quotient, is order-preserving. 


The content of this chapter is unlike anything that has been covered in the first part. In this chapter we try to answer the most significant difference between the paths of the first part and the second-order paths introduced in this second part, namely the need to work with class representatives in second-order paths. 

The $\Sigma^{\boldsymbol{\mathcal{A}}^{(2)}}$-congruence we present in this chapter aims to avoid this forced choice. However, these conditions materialise in a minimal congruence that contains the equations we need. Otherwise, we might be collapsing too many second-order paths and lose control over the equations. This new relation will help, for example, to solve the associativity problem for the $0$-composition. See Remark~\ref{RDCompZ}.

We will now define the relation we will be working with in this chapter.

\begin{restatable}{definition}{DDUps}
\label{DDUps}
\index{Upsilon!second-order!$\Upsilon^{(1)}$}
 We define $\Upsilon^{(1)}=(\Upsilon^{(1)}_{s})_{s\in S}$ to be the relation on $\mathrm{Pth}_{\boldsymbol{\mathcal{A}}^{(2)}}$ consisting exactly of the following pairs of second-order paths
\begin{enumerate}
\item For every sort $s\in S$ and every second-order path $\mathfrak{P}^{(2)}$ in $\mathrm{Pth}_{\boldsymbol{\mathcal{A}}^{(2)},s}$, 
\[
\left(
\mathfrak{P}^{(2)},
\mathfrak{P}^{(2)}\circ^{0\mathbf{Pth}_{\boldsymbol{\mathcal{A}}}^{(2)}}_{s}
\mathrm{sc}^{0\mathbf{Pth}_{\boldsymbol{\mathcal{A}}^{(2)}}}_{s}\left(
\mathfrak{P}^{(2)}
\right)
\right)
\in\Upsilon^{(1)}_{s};
\]
\item For every sort $s\in S$ and every second-order path $\mathfrak{P}^{(2)}$ in $\mathrm{Pth}_{\boldsymbol{\mathcal{A}}^{(2)},s}$, 
\[
\left(
\mathfrak{P}^{(2)},
\mathrm{tg}^{0\mathbf{Pth}_{\boldsymbol{\mathcal{A}}^{(2)}}}_{s}\left(
\mathfrak{P}^{(2)}
\right)\circ^{0\mathbf{Pth}_{\boldsymbol{\mathcal{A}}}^{(2)}}_{s}
\mathfrak{P}^{(2)}
\right)
\in\Upsilon^{(1)}_{s};
\]
\item For every sort $s\in S$ and every second-order paths $\mathfrak{P}^{(2)}$, $\mathfrak{Q}^{(2)}$ and $\mathfrak{R}^{(2)}$ in $\mathrm{Pth}_{\boldsymbol{\mathcal{A}}^{(2)},s}$ with
\allowdisplaybreaks
\begin{align*}
\mathrm{sc}^{(0,2)}_{s}\left(
\mathfrak{R}^{(2)}
\right)
&=
\mathrm{tg}^{(0,2)}_{s}\left(
\mathfrak{Q}^{(2)}
\right);
&
\mathrm{sc}^{(0,2)}_{s}\left(
\mathfrak{Q}^{(2)}
\right)
&=
\mathrm{tg}^{(0,2)}_{s}\left(
\mathfrak{R}^{(2)}
\right);
\end{align*}
then
\allowdisplaybreaks
\begin{multline*}
\left(
\mathfrak{R}^{(2)}
\circ^{0\mathbf{Pth}_{\boldsymbol{\mathcal{A}}^{(2)}}}_{s}
\left(
\mathfrak{Q}^{(2)}
\circ^{0\mathbf{Pth}_{\boldsymbol{\mathcal{A}}^{(2)}}}_{s}
\mathfrak{P}^{(2)}
\right)
,
\right.
\\
\left.
\left(
\mathfrak{R}^{(2)}
\circ^{0\mathbf{Pth}_{\boldsymbol{\mathcal{A}}^{(2)}}}_{s}
\mathfrak{Q}^{(2)}
\right)
\circ^{0\mathbf{Pth}_{\boldsymbol{\mathcal{A}}^{(2)}}}_{s}
\mathfrak{P}^{(2)}
\right)
\in\Upsilon^{(1)}_{s};
\end{multline*}
\item For every word $\mathbf{s}\in S^{\star}$ and $s\in S$, for every operation symbol $\sigma\in \Sigma_{\mathbf{s},s}$, for every two families of second-order paths $(\mathfrak{P}^{(2)}_{j})_{j\in\bb{\mathbf{s}}}$ and $(\mathfrak{Q}^{(2)}_{j})_{j\in\bb{\mathbf{s}}}$  in $\mathrm{Pth}_{\boldsymbol{\mathcal{A}}^{(2)},\mathbf{s}}$ satisfying that, for every $j\in\bb{\mathbf{s}}$,
\[
\mathrm{sc}^{(0,2)}_{s_{j}}\left(
\mathfrak{Q}^{(2)}_{j}
\right)
=
\mathrm{tg}^{(0,2)}_{s_{j}}\left(
\mathfrak{P}^{(2)}_{j}
\right),
\]
then
\allowdisplaybreaks
\begin{multline*}
\left(
\sigma^{\mathbf{Pth}_{\boldsymbol{\mathcal{A}}^{(2)}}}\left(
\left(
\mathfrak{Q}^{(2)}_{j}
\circ^{0\mathbf{Pth}_{\boldsymbol{\mathcal{A}}^{(2)}}}_{s_{j}}
\mathfrak{P}^{(2)}_{j}
\right)_{j\in\bb{\mathbf{s}}}
\right),
\right.
\\
\left.
\sigma^{\mathbf{Pth}_{\boldsymbol{\mathcal{A}}^{(2)}}}\left(
\left(
\mathfrak{Q}^{(2)}_{j}
\right)_{j\in\bb{\mathbf{s}}}
\right)
\circ^{0\mathbf{Pth}_{\boldsymbol{\mathcal{A}}^{(2)}}}_{s}
\sigma^{\mathbf{Pth}_{\boldsymbol{\mathcal{A}}^{(2)}}}\left(
\left(
\mathfrak{P}^{(2)}_{j}
\right)_{j\in\bb{\mathbf{s}}}
\right)
\right)
\in\Upsilon^{(1)}_{s}.
\end{multline*}
\end{enumerate}
This completes the definition of $\Upsilon^{(1)}$.
\end{restatable}

We can draw some important consequences directly from the above definition. 
The first of these consequences would be that, if two second-order paths are related by $\Upsilon^{(1)}$, then they have the same length, $([1],2)$-source and $([1],2)$-target.

\begin{restatable}{lemma}{LDUps}
\label{LDUps} Let $s$ be a sort in $S$ and $(\mathfrak{P}^{(2)},\mathfrak{Q}^{(2)})\in\Upsilon^{(1)}_{s}$. The following statements hold.
\begin{enumerate}
\item[(i)] $\bb{\mathfrak{P}^{(2)}}=\bb{\mathfrak{Q}^{(2)}}$;
\item[(ii)] $\mathrm{sc}_{s}^{([1],2)}(\mathfrak{P}^{(2)})=\mathrm{sc}_{s}^{([1],2)}(\mathfrak{Q}^{(2)})$;
\item[(iii)] $\mathrm{tg}_{s}^{([1],2)}(\mathfrak{P}^{(2)})=\mathrm{tg}_{s}^{([1],2)}(\mathfrak{Q}^{(2)})$.
\end{enumerate}
\end{restatable}
\begin{proof}
The proof involves checking each case used in the definition of $\Upsilon^{(1)}$.

Regarding item (1), we have that, for every sort $s\in S$ and every second-order path $\mathfrak{P}^{(2)}$ in $\mathrm{Pth}_{\boldsymbol{\mathcal{A}}^{(2)},s}$, 
\[
\left(
\mathfrak{P}^{(2)},
\mathfrak{P}^{(2)}\circ^{0\mathbf{Pth}_{\boldsymbol{\mathcal{A}}}^{(2)}}_{s}
\mathrm{sc}^{0\mathbf{Pth}_{\boldsymbol{\mathcal{A}}^{(2)}}}_{s}\left(
\mathfrak{P}^{(2)}
\right)
\right)
\in\Upsilon^{(1)}_{s};
\]

Note that the following chain of equalities hold
\begin{align*}
\left\lvert
\mathfrak{P}^{(2)}\circ^{0\mathbf{Pth}_{\boldsymbol{\mathcal{A}}}^{(2)}}_{s}
\mathrm{sc}^{0\mathbf{Pth}_{\boldsymbol{\mathcal{A}}^{(2)}}}_{s}\left(
\mathfrak{P}^{(2)}
\right)
\right\rvert
&=
\left\lvert
\mathfrak{P}^{(2)}
\right\rvert
+
\left\lvert
\mathrm{sc}^{0\mathbf{Pth}_{\boldsymbol{\mathcal{A}}^{(2)}}}_{s}\left(
\mathfrak{P}^{(2)}
\right)
\right\rvert
\tag{1}
\\&=
\left\lvert
\mathfrak{P}^{(2)}
\right\rvert.
\tag{2}
\end{align*}

In the just stated chain of equalities, the first equality follows from Claim~\ref{CDPthCatAlgCompZ}; the second equality follows from Claim~\ref{CDPthCatAlgScZ}.

Moreover, the following chain of equalities hold
\begin{flushleft}
$\mathrm{sc}^{([1],2)}_{s}\left(
\mathfrak{P}^{(2)}\circ^{0\mathbf{Pth}_{\boldsymbol{\mathcal{A}}}^{(2)}}_{s}
\mathrm{sc}^{0\mathbf{Pth}_{\boldsymbol{\mathcal{A}}^{(2)}}}_{s}\left(
\mathfrak{P}^{(2)}
\right)\right)$
\allowdisplaybreaks
\begin{align*}
\qquad
&=
\mathrm{sc}^{([1],2)}_{s}\left(
\mathfrak{P}^{(2)}
\right)
\circ^{0[\mathbf{PT}_{\boldsymbol{\mathcal{A}}}]}_{s}
\mathrm{sc}^{([1],2)}_{s}\left(
\mathrm{sc}^{0\mathbf{Pth}_{\boldsymbol{\mathcal{A}}^{(2)}}}_{s}\left(
\mathfrak{P}^{(2)}
\right)
\right)
\tag{1}
\\&=
\mathrm{sc}^{([1],2)}_{s}\left(
\mathfrak{P}^{(2)}
\right)
\circ^{0[\mathbf{PT}_{\boldsymbol{\mathcal{A}}}]}_{s}
\mathrm{sc}^{([1],2)}_{s}\left(
\mathrm{ip}^{(2,[1])\sharp}_{s}\left(
\mathrm{sc}^{0[\mathbf{PT}_{\boldsymbol{\mathcal{A}}}]}_{s}\left(
\mathrm{sc}^{([1],2)}_{s}\left(
\mathfrak{P}^{(2)}
\right)
\right)
\right)
\right)
\tag{2}
\\&=
\mathrm{sc}^{([1],2)}_{s}\left(
\mathfrak{P}^{(2)}
\right)
\circ^{0[\mathbf{PT}_{\boldsymbol{\mathcal{A}}}]}_{s}
\mathrm{sc}^{0[\mathbf{PT}_{\boldsymbol{\mathcal{A}}}]}_{s}\left(
\mathrm{sc}^{([1],2)}_{s}\left(
\mathfrak{P}^{(2)}
\right)
\right)
\tag{3}
\\&=
\mathrm{sc}^{([1],2)}_{s}\left(
\mathfrak{P}^{(2)}
\right).
\tag{4}
\end{align*}
\end{flushleft}

In the just stated chain of equalities, the first equality follows from Claim~\ref{CDPthCatAlgCompZ}; the second equality follows from Claim~\ref{CDPthCatAlgScZ}; the third equality follows from Proposition~\ref{PDBasicEq}; finally, the last equality follows from Proposition~\ref{PPTQVarA5}.

Now, regarding the $([1],2)$-target, the following chain of equalities hold
\begin{flushleft}
$\mathrm{tg}^{([1],2)}_{s}\left(
\mathfrak{P}^{(2)}\circ^{0\mathbf{Pth}_{\boldsymbol{\mathcal{A}}}^{(2)}}_{s}
\mathrm{sc}^{0\mathbf{Pth}_{\boldsymbol{\mathcal{A}}^{(2)}}}_{s}\left(
\mathfrak{P}^{(2)}
\right)\right)$
\allowdisplaybreaks
\begin{align*}
\qquad
&=
\mathrm{tg}^{([1],2)}_{s}\left(
\mathfrak{P}^{(2)}
\right)
\circ^{0[\mathbf{PT}_{\boldsymbol{\mathcal{A}}}]}_{s}
\mathrm{tg}^{([1],2)}_{s}\left(
\mathrm{sc}^{0\mathbf{Pth}_{\boldsymbol{\mathcal{A}}^{(2)}}}_{s}\left(
\mathfrak{P}^{(2)}
\right)
\right)
\tag{1}
\\&=
\mathrm{tg}^{([1],2)}_{s}\left(
\mathfrak{P}^{(2)}
\right)
\circ^{0[\mathbf{PT}_{\boldsymbol{\mathcal{A}}}]}_{s}
\mathrm{tg}^{([1],2)}_{s}\left(
\mathrm{ip}^{(2,[1])\sharp}_{s}\left(
\mathrm{sc}^{0[\mathbf{PT}_{\boldsymbol{\mathcal{A}}}]}_{s}\left(
\mathrm{tg}^{([1],2)}_{s}\left(
\mathfrak{P}^{(2)}
\right)
\right)
\right)
\right)
\tag{2}
\\&=
\mathrm{tg}^{([1],2)}_{s}\left(
\mathfrak{P}^{(2)}
\right)
\circ^{0[\mathbf{PT}_{\boldsymbol{\mathcal{A}}}]}_{s}
\mathrm{sc}^{0[\mathbf{PT}_{\boldsymbol{\mathcal{A}}}]}_{s}\left(
\mathrm{tg}^{([1],2)}_{s}\left(
\mathfrak{P}^{(2)}
\right)
\right)
\tag{3}
\\&=
\mathrm{tg}^{([1],2)}_{s}\left(
\mathfrak{P}^{(2)}
\right).
\tag{4}
\end{align*}
\end{flushleft}

In the just stated chain of equalities, the first equality follows from Claim~\ref{CDPthCatAlgCompZ}; the second equality follows from Claim~\ref{CDPthCatAlgScZ}; the third equality follows from Proposition~\ref{PDBasicEq}; finally, the last equality follows from Proposition~\ref{PPTQVarA5}.

This completes the proof of the first item.

Regarding item (2), we have that, for every sort $s\in S$ and every second-order path $\mathfrak{P}^{(2)}$ in $\mathrm{Pth}_{\boldsymbol{\mathcal{A}}^{(2)},s}$, 
\[
\left(
\mathfrak{P}^{(2)},
\mathrm{tg}^{0\mathbf{Pth}_{\boldsymbol{\mathcal{A}}^{(2)}}}_{s}\left(
\mathfrak{P}^{(2)}
\right)\circ^{0\mathbf{Pth}_{\boldsymbol{\mathcal{A}}}^{(2)}}_{s}
\mathfrak{P}^{(2)}
\right)
\in\Upsilon^{(1)}_{s}.
\]

This case follows by a similar argument to that in the proof presented above for the first case. 

This completes the proof of the second item.

Regarding item (3), we have that, for every sort $s\in S$ and every second-order paths $\mathfrak{P}^{(2)}$, $\mathfrak{Q}^{(2)}$ and $\mathfrak{R}^{(2)}$ in $\mathrm{Pth}_{\boldsymbol{\mathcal{A}}^{(2)},s}$ with
\allowdisplaybreaks
\begin{align*}
\mathrm{sc}^{(0,2)}_{s}\left(
\mathfrak{R}^{(2)}
\right)
&=
\mathrm{tg}^{(0,2)}_{s}\left(
\mathfrak{Q}^{(2)}
\right);
&
\mathrm{sc}^{(0,2)}_{s}\left(
\mathfrak{Q}^{(2)}
\right)
&=
\mathrm{tg}^{(0,2)}_{s}\left(
\mathfrak{R}^{(2)}
\right);
\end{align*}
then the following pair is in $\Upsilon^{(1)}_{s}$,
\allowdisplaybreaks
\begin{multline*}
\left(
\mathfrak{R}^{(2)}
\circ^{0\mathbf{Pth}_{\boldsymbol{\mathcal{A}}^{(2)}}}_{s}
\left(
\mathfrak{Q}^{(2)}
\circ^{0\mathbf{Pth}_{\boldsymbol{\mathcal{A}}^{(2)}}}_{s}
\mathfrak{P}^{(2)}
\right)
,
\right.
\\
\left.
\left(
\mathfrak{R}^{(2)}
\circ^{0\mathbf{Pth}_{\boldsymbol{\mathcal{A}}^{(2)}}}_{s}
\mathfrak{Q}^{(2)}
\right)
\circ^{0\mathbf{Pth}_{\boldsymbol{\mathcal{A}}^{(2)}}}_{s}
\mathfrak{P}^{(2)}
\right)
\in\Upsilon^{(1)}_{s};
\end{multline*}

Note that, in virtue of Claim~\ref{CDPthCatAlgCompZ}, the following chain of equalities hold
\allowdisplaybreaks
\begin{multline*}
\left\lvert
\mathfrak{R}^{(2)}
\circ^{0\mathbf{Pth}_{\boldsymbol{\mathcal{A}}^{(2)}}}_{s}
\left(
\mathfrak{Q}^{(2)}
\circ^{0\mathbf{Pth}_{\boldsymbol{\mathcal{A}}^{(2)}}}_{s}
\mathfrak{P}^{(2)}
\right)
\right\rvert
\\=
\left\lvert
\mathfrak{R}^{(2)}
\right\rvert
+
\left\lvert
\mathfrak{Q}^{(2)}
\right\rvert
+
\left\lvert
\mathfrak{P}^{(2)}
\right\rvert
=
\left\lvert
\left(
\mathfrak{R}^{(2)}
\circ^{0\mathbf{Pth}_{\boldsymbol{\mathcal{A}}^{(2)}}}_{s}
\mathfrak{Q}^{(2)}
\right)
\circ^{0\mathbf{Pth}_{\boldsymbol{\mathcal{A}}^{(2)}}}_{s}
\mathfrak{P}^{(2)}
\right\rvert.
\end{multline*}

Regarding the $([1],2)$-source, the following chain of equalities hold
\begin{flushleft}
$\mathrm{sc}^{([1],2)}_{s}\left(
\mathfrak{R}^{(2)}
\circ^{0\mathbf{Pth}_{\boldsymbol{\mathcal{A}}^{(2)}}}_{s}
\left(
\mathfrak{Q}^{(2)}
\circ^{0\mathbf{Pth}_{\boldsymbol{\mathcal{A}}^{(2)}}}_{s}
\mathfrak{P}^{(2)}
\right)
\right)$
\allowdisplaybreaks
\begin{align*}
\qquad&=
\mathrm{sc}^{([1],2)}_{s}\left(
\mathfrak{R}^{(2)}
\right)
\circ^{0[\mathbf{PT}_{\boldsymbol{\mathcal{A}}}]}_{s}
\left(
\mathrm{sc}^{([1],2)}_{s}\left(
\mathfrak{Q}^{(2)}
\right)
\circ^{0[\mathbf{PT}_{\boldsymbol{\mathcal{A}}}]}_{s}
\mathrm{sc}^{([1],2)}_{s}\left(
\mathfrak{P}^{(2)}
\right)
\right)
\tag{1}
\\&=
\left(
\mathrm{sc}^{([1],2)}_{s}\left(
\mathfrak{R}^{(2)}
\right)
\circ^{0[\mathbf{PT}_{\boldsymbol{\mathcal{A}}}]}_{s}
\mathrm{sc}^{([1],2)}_{s}\left(
\mathfrak{Q}^{(2)}
\right)
\right)
\circ^{0[\mathbf{PT}_{\boldsymbol{\mathcal{A}}}]}_{s}
\mathrm{sc}^{([1],2)}_{s}\left(
\mathfrak{P}^{(2)}
\right)
\tag{2}
\\&=
\mathrm{sc}^{([1],2)}_{s}\left(
\left(
\mathfrak{R}^{(2)}
\circ^{0\mathbf{Pth}_{\boldsymbol{\mathcal{A}}^{(2)}}}_{s}
\mathfrak{Q}^{(2)}
\right)
\circ^{0\mathbf{Pth}_{\boldsymbol{\mathcal{A}}^{(2)}}}_{s}
\mathfrak{P}^{(2)}
\right).
\tag{3}
\end{align*}
\end{flushleft}

In the just stated chain of equalities, the first equality follows from Claim~\ref{CDPthCatAlgCompZ}; the second equality follows from Proposition~\ref{PPTQVarA6}; finally, the last equality follows from Claim~\ref{CDPthCatAlgCompZ}.

Regarding the $([1],2)$-target, the following chain of equalities hold
\begin{flushleft}
$\mathrm{tg}^{([1],2)}_{s}\left(
\mathfrak{R}^{(2)}
\circ^{0\mathbf{Pth}_{\boldsymbol{\mathcal{A}}^{(2)}}}_{s}
\left(
\mathfrak{Q}^{(2)}
\circ^{0\mathbf{Pth}_{\boldsymbol{\mathcal{A}}^{(2)}}}_{s}
\mathfrak{P}^{(2)}
\right)
\right)$
\allowdisplaybreaks
\begin{align*}
\qquad&=
\mathrm{tg}^{([1],2)}_{s}\left(
\mathfrak{R}^{(2)}
\right)
\circ^{0[\mathbf{PT}_{\boldsymbol{\mathcal{A}}}]}_{s}
\left(
\mathrm{tg}^{([1],2)}_{s}\left(
\mathfrak{Q}^{(2)}
\right)
\circ^{0[\mathbf{PT}_{\boldsymbol{\mathcal{A}}}]}_{s}
\mathrm{tg}^{([1],2)}_{s}\left(
\mathfrak{P}^{(2)}
\right)
\right)
\tag{1}
\\&=
\left(
\mathrm{tg}^{([1],2)}_{s}\left(
\mathfrak{R}^{(2)}
\right)
\circ^{0[\mathbf{PT}_{\boldsymbol{\mathcal{A}}}]}_{s}
\mathrm{tg}^{([1],2)}_{s}\left(
\mathfrak{Q}^{(2)}
\right)
\right)
\circ^{0[\mathbf{PT}_{\boldsymbol{\mathcal{A}}}]}_{s}
\mathrm{tg}^{([1],2)}_{s}\left(
\mathfrak{P}^{(2)}
\right)
\tag{2}
\\&=
\mathrm{tg}^{([1],2)}_{s}\left(
\left(
\mathfrak{R}^{(2)}
\circ^{0\mathbf{Pth}_{\boldsymbol{\mathcal{A}}^{(2)}}}_{s}
\mathfrak{Q}^{(2)}
\right)
\circ^{0\mathbf{Pth}_{\boldsymbol{\mathcal{A}}^{(2)}}}_{s}
\mathfrak{P}^{(2)}
\right).
\tag{3}
\end{align*}
\end{flushleft}

In the just stated chain of equalities, the first equality follows from Claim~\ref{CDPthCatAlgCompZ}; the second equality follows from Proposition~\ref{PPTQVarA6}; finally, the last equality follows from Claim~\ref{CDPthCatAlgCompZ}.

This completes the proof of the third item.

Regarding item (4), we have that, for every word $\mathbf{s}\in S^{\star}$ and $s\in S$, for every operation symbol $\sigma\in \Sigma_{\mathbf{s},s}$, for every families of second-order paths $(\mathfrak{P}^{(2)}_{j})_{j\in\bb{\mathbf{s}}}$ and $(\mathfrak{Q}^{(2)}_{j})_{j\in\bb{\mathbf{s}}}$  in $\mathrm{Pth}_{\boldsymbol{\mathcal{A}}^{(2)},\mathbf{s}}$ satisfying that, for every $j\in\bb{\mathbf{s}}$,
\[
\mathrm{sc}^{(0,2)}_{s_{j}}\left(
\mathfrak{Q}^{(2)}_{j}
\right)
=
\mathrm{tg}^{(0,2)}_{s_{j}}\left(
\mathfrak{P}^{(2)}_{j}
\right),
\]
then the following pair is in $\Upsilon^{(1)}_{s}$, 
\allowdisplaybreaks
\begin{multline*}
\left(
\sigma^{\mathbf{Pth}_{\boldsymbol{\mathcal{A}}^{(2)}}}\left(
\left(
\mathfrak{Q}^{(2)}_{j}
\circ^{0\mathbf{Pth}_{\boldsymbol{\mathcal{A}}^{(2)}}}_{s_{j}}
\mathfrak{P}^{(2)}_{j}
\right)_{j\in\bb{\mathbf{s}}}
\right),
\right.
\\
\left.
\sigma^{\mathbf{Pth}_{\boldsymbol{\mathcal{A}}^{(2)}}}\left(
\left(
\mathfrak{Q}^{(2)}_{j}
\right)_{j\in\bb{\mathbf{s}}}
\right)
\circ^{0\mathbf{Pth}_{\boldsymbol{\mathcal{A}}^{(2)}}}_{s}
\sigma^{\mathbf{Pth}_{\boldsymbol{\mathcal{A}}^{(2)}}}\left(
\left(
\mathfrak{P}^{(2)}_{j}
\right)_{j\in\bb{\mathbf{s}}}
\right)
\right)
\in\Upsilon^{(1)}_{s}.
\end{multline*}

Note that the following chain of equalities hold
\begin{flushleft}
$\left\lvert
\sigma^{\mathbf{Pth}_{\boldsymbol{\mathcal{A}}^{(2)}}}\left(
\left(
\mathfrak{Q}^{(2)}_{j}
\circ^{0\mathbf{Pth}_{\boldsymbol{\mathcal{A}}^{(2)}}}_{s_{j}}
\mathfrak{P}^{(2)}_{j}
\right)_{j\in\bb{\mathbf{s}}}
\right)
\right\rvert $
\allowdisplaybreaks
\begin{align*}
\qquad
&=
\sum_{j\in\bb{\mathbf{s}}}\left\lvert
\mathfrak{Q}^{(2)}_{j}
\circ^{0\mathbf{Pth}_{\boldsymbol{\mathcal{A}}^{(2)}}}_{s_{j}}
\mathfrak{P}^{(2)}_{j}
\right\rvert
\tag{1}
\\&=
\sum_{j\in\bb{\mathbf{s}}}
\left(\left\lvert
\mathfrak{Q}^{(2)}_{j}
\right\rvert
+
\left\lvert
\mathfrak{P}^{(2)}_{j}
\right\rvert
\right)
\tag{2}
\\&=
\sum_{j\in\bb{\mathbf{s}}}\left\lvert
\mathfrak{Q}^{(2)}_{j}
\right\rvert
+
\sum_{j\in\bb{\mathbf{s}}}
\left\lvert
\mathfrak{P}^{(2)}_{j}
\right\rvert
\tag{3}
\\&=
\left\lvert
\sigma^{\mathbf{Pth}_{\boldsymbol{\mathcal{A}}^{(2)}}}\left(
\left(
\mathfrak{Q}^{(2)}_{j}
\right)_{j\in\bb{\mathbf{s}}}
\right)
\right\rvert
+
\left\lvert
\sigma^{\mathbf{Pth}_{\boldsymbol{\mathcal{A}}^{(2)}}}\left(
\left(
\mathfrak{P}^{(2)}_{j}
\right)_{j\in\bb{\mathbf{s}}}
\right)
\right\rvert
\tag{4}
\\&=
\left\lvert
\sigma^{\mathbf{Pth}_{\boldsymbol{\mathcal{A}}^{(2)}}}\left(
\left(
\mathfrak{Q}^{(2)}_{j}
\right)_{j\in\bb{\mathbf{s}}}
\right)
\circ^{0\mathbf{Pth}_{\boldsymbol{\mathcal{A}}^{(2)}}}_{s}
\sigma^{\mathbf{Pth}_{\boldsymbol{\mathcal{A}}^{(2)}}}\left(
\left(
\mathfrak{P}^{(2)}_{j}
\right)_{j\in\bb{\mathbf{s}}}
\right)
\right\rvert.
\tag{5}
\end{align*}
\end{flushleft}

In the just stated chain of equalities, the first equality follows in virtue of Claim~\ref{CDPthSigma}, the second equality follows from Claim~\ref{CDPthCatAlgCompZ}; the third equality follows from the associativity and commutativity of the sum; the fourth equality follows from Claim~\ref{CDPthSigma}; finally, the last equality follows from Claim~\ref{CDPthCatAlgCompZ}.

Regarding the $([1],2)$-source, the following chain of equalities hold
\begin{flushleft}
$\mathrm{sc}^{([1],2)}_{s}\left(
\sigma^{\mathbf{Pth}_{\boldsymbol{\mathcal{A}}^{(2)}}}\left(
\left(
\mathfrak{Q}^{(2)}_{j}
\circ^{0\mathbf{Pth}_{\boldsymbol{\mathcal{A}}^{(2)}}}_{s_{j}}
\mathfrak{P}^{(2)}_{j}
\right)_{j\in\bb{\mathbf{s}}}
\right)
\right)$
\allowdisplaybreaks
\begin{align*}
&=
\sigma^{[\mathbf{PT}_{\boldsymbol{\mathcal{A}}}]}\left(
\left(
\mathrm{sc}^{([1],2)}_{s_{j}}\left(
\mathfrak{Q}^{(2)}_{j}
\circ^{0\mathbf{Pth}_{\boldsymbol{\mathcal{A}}^{(2)}}}_{s_{j}}
\mathfrak{P}^{(2)}_{j}
\right)
\right)_{j\in\bb{\mathbf{s}}}
\right)
\tag{1}
\\&=
\sigma^{[\mathbf{PT}_{\boldsymbol{\mathcal{A}}}]}\left(
\left(
\mathrm{sc}^{([1],2)}_{s_{j}}\left(
\mathfrak{Q}^{(2)}_{j}
\right)
\circ^{0[\mathbf{PT}_{\boldsymbol{\mathcal{A}}}]}_{s_{j}}
\mathrm{sc}^{([1],2)}_{s_{j}}\left(
\mathfrak{P}^{(2)}_{j}
\right)
\right)_{j\in\bb{\mathbf{s}}}
\right)
\tag{2}
\\&=
\sigma^{[\mathbf{PT}_{\boldsymbol{\mathcal{A}}}]}\left(
\left(
\mathrm{sc}^{([1],2)}_{s_{j}}\left(
\mathfrak{Q}^{(2)}_{j}
\right)
\right)_{j\in\bb{\mathbf{s}}}
\right)
\circ^{0[\mathbf{PT}_{\boldsymbol{\mathcal{A}}}]}_{s}
\sigma^{[\mathbf{PT}_{\boldsymbol{\mathcal{A}}}]}\left(
\left(
\mathrm{sc}^{([1],2)}_{s_{j}}\left(
\mathfrak{P}^{(2)}_{j}
\right)
\right)_{j\in\bb{\mathbf{s}}}
\right)
\tag{3}
\\&=
\mathrm{sc}^{([1],2)}_{s}\left(
\sigma^{[\mathbf{PT}_{\boldsymbol{\mathcal{A}}}]}\left(
\left(
\mathfrak{Q}^{(2)}_{j}
\right)_{j\in\bb{\mathbf{s}}}
\right)
\right)
\circ^{0[\mathbf{PT}_{\boldsymbol{\mathcal{A}}}]}_{s}
\mathrm{sc}^{([1],2)}_{s}\left(
\sigma^{[\mathbf{PT}_{\boldsymbol{\mathcal{A}}}]}\left(
\left(
\mathfrak{P}^{(2)}_{j}
\right)_{j\in\bb{\mathbf{s}}}
\right)
\right)
\tag{4}
\\&=
\mathrm{sc}^{([1],2)}_{s}\left(
\sigma^{\mathbf{Pth}_{\boldsymbol{\mathcal{A}}^{(2)}}}\left(
\left(
\mathfrak{Q}^{(2)}_{j}
\right)_{j\in\bb{\mathbf{s}}}
\right)
\circ^{0\mathbf{Pth}_{\boldsymbol{\mathcal{A}}^{(2)}}}_{s}
\sigma^{\mathbf{Pth}_{\boldsymbol{\mathcal{A}}^{(2)}}}\left(
\left(
\mathfrak{P}^{(2)}_{j}
\right)_{j\in\bb{\mathbf{s}}}
\right)
\right).
\tag{5}
\end{align*}
\end{flushleft}

In the just stated chain of equalities, the first equality follows from Claim~\ref{CDPthSigma}; the second equality follows from Claim~\ref{CDPthCatAlgCompZ}; the third equality follows from Proposition~\ref{PPTQVarA8}; the fourth equality follows from Claim~\ref{CDPthSigma}; finally, the last equality follows from Claim~\ref{CDPthCatAlgCompZ}.

Regarding the $([1],2)$-target, the following chain of equalities hold
\begin{flushleft}
$\mathrm{tg}^{([1],2)}_{s}\left(
\sigma^{\mathbf{Pth}_{\boldsymbol{\mathcal{A}}^{(2)}}}\left(
\left(
\mathfrak{Q}^{(2)}_{j}
\circ^{0\mathbf{Pth}_{\boldsymbol{\mathcal{A}}^{(2)}}}_{s_{j}}
\mathfrak{P}^{(2)}_{j}
\right)_{j\in\bb{\mathbf{s}}}
\right)
\right)$
\allowdisplaybreaks
\begin{align*}
&=
\sigma^{[\mathbf{PT}_{\boldsymbol{\mathcal{A}}}]}\left(
\left(
\mathrm{tg}^{([1],2)}_{s_{j}}\left(
\mathfrak{Q}^{(2)}_{j}
\circ^{0\mathbf{Pth}_{\boldsymbol{\mathcal{A}}^{(2)}}}_{s_{j}}
\mathfrak{P}^{(2)}_{j}
\right)
\right)_{j\in\bb{\mathbf{s}}}
\right)
\tag{1}
\\&=
\sigma^{[\mathbf{PT}_{\boldsymbol{\mathcal{A}}}]}\left(
\left(
\mathrm{tg}^{([1],2)}_{s_{j}}\left(
\mathfrak{Q}^{(2)}_{j}
\right)
\circ^{0[\mathbf{PT}_{\boldsymbol{\mathcal{A}}}]}_{s_{j}}
\mathrm{tg}^{([1],2)}_{s_{j}}\left(
\mathfrak{P}^{(2)}_{j}
\right)
\right)_{j\in\bb{\mathbf{s}}}
\right)
\tag{2}
\\&=
\sigma^{[\mathbf{PT}_{\boldsymbol{\mathcal{A}}}]}\left(
\left(
\mathrm{tg}^{([1],2)}_{s_{j}}\left(
\mathfrak{Q}^{(2)}_{j}
\right)
\right)_{j\in\bb{\mathbf{s}}}
\right)
\circ^{0[\mathbf{PT}_{\boldsymbol{\mathcal{A}}}]}_{s}
\sigma^{[\mathbf{PT}_{\boldsymbol{\mathcal{A}}}]}\left(
\left(
\mathrm{tg}^{([1],2)}_{s_{j}}\left(
\mathfrak{P}^{(2)}_{j}
\right)
\right)_{j\in\bb{\mathbf{s}}}
\right)
\tag{3}
\\&=
\mathrm{tg}^{([1],2)}_{s}\left(
\sigma^{[\mathbf{PT}_{\boldsymbol{\mathcal{A}}}]}\left(
\left(
\mathfrak{Q}^{(2)}_{j}
\right)_{j\in\bb{\mathbf{s}}}
\right)
\right)
\circ^{0[\mathbf{PT}_{\boldsymbol{\mathcal{A}}}]}_{s}
\mathrm{tg}^{([1],2)}_{s}\left(
\sigma^{[\mathbf{PT}_{\boldsymbol{\mathcal{A}}}]}\left(
\left(
\mathfrak{P}^{(2)}_{j}
\right)_{j\in\bb{\mathbf{s}}}
\right)
\right)
\tag{4}
\\&=
\mathrm{tg}^{([1],2)}_{s}\left(
\sigma^{\mathbf{Pth}_{\boldsymbol{\mathcal{A}}^{(2)}}}\left(
\left(
\mathfrak{Q}^{(2)}_{j}
\right)_{j\in\bb{\mathbf{s}}}
\right)
\circ^{0\mathbf{Pth}_{\boldsymbol{\mathcal{A}}^{(2)}}}_{s}
\sigma^{\mathbf{Pth}_{\boldsymbol{\mathcal{A}}^{(2)}}}\left(
\left(
\mathfrak{P}^{(2)}_{j}
\right)_{j\in\bb{\mathbf{s}}}
\right)
\right).
\tag{5}
\end{align*}
\end{flushleft}

In the just stated chain of equalities, the first equality follows from Claim~\ref{CDPthSigma}; the second equality follows from Claim~\ref{CDPthCatAlgCompZ}; the third equality follows from Proposition~\ref{PPTQVarA8}; the fourth equality follows from Claim~\ref{CDPthSigma}; finally, the last equality follows from Claim~\ref{CDPthCatAlgCompZ}.

This completes the proof of the last item. 

This completes the proof.
\end{proof}

Next we see that, in the case of dealing with $(2,[1])$-identity second-order paths, the  relation $\Upsilon^{(1)}$ is the diagonal relation.

\begin{restatable}{corollary}{CDUps}
\label{CDUps} Let $s$ be a sort in $S$ and $(\mathfrak{P}^{(2)},\mathfrak{Q}^{(2)})\in\Upsilon^{(1)}_{s}$. If $\mathfrak{P}^{(2)}$ or $\mathfrak{Q}^{(2)}$ is a $(2,[1])$-identity second-order path, then $\mathfrak{P}^{(2)}=\mathfrak{Q}^{(2)}$.
\end{restatable}
\begin{proof}
It follows from Propositions~\ref{PDUIpVarA5},~\ref{PDUIpVarA6} and~\ref{PDUIpVarA8}.
\end{proof}

\section{
\texorpdfstring
{The congruence $\Upsilon^{[1]}$ on $\mathbf{Pth}_{\boldsymbol{\mathcal{A}}^{(2)}}$}
{The Upsilon congruence}
}

We next show that the properties that hold for the relation $\Upsilon^{(1)}$, set up in Definition~\ref{DDUps}, also hold for the smallest $\Sigma^{\boldsymbol{\mathcal{A}}^{(2)}}$-congruence on $\mathbf{Pth}_{\boldsymbol{\mathcal{A}}^{(2)}}$ containing $\Upsilon^{(1)}$. In this regard, we recover the notation introduced in Definition~\ref{DCongOpInt}.

\begin{restatable}{definition}{DDUpsCong}
\index{Upsilon!second-order!$\Upsilon^{[1]}$}
\label{DDUpsCong} We denote by $\Upsilon^{[1]}$ the smallest $\Sigma^{\boldsymbol{\mathcal{A}}^{(2)}}$-congruence on $\mathbf{Pth}_{\boldsymbol{\mathcal{A}}^{(2)}}$ containing $\Upsilon^{(1)}$, i.e., 
\[
\Upsilon^{[1]}=\mathrm{Cg}_{\mathbf{Pth}_{\boldsymbol{\mathcal{A}}^{(2)}}}\left(
\Upsilon^{(1)}
\right).
\]
\end{restatable}

\begin{remark}\label{RDUpsCong}
Let $\mathrm{C}_{\mathbf{Pth}_{\boldsymbol{\mathcal{A}}^{(2)}}}$ be the operator on $\mathrm{Pth}_{\boldsymbol{\mathcal{A}}^{(2)}}\times \mathrm{Pth}_{\boldsymbol{\mathcal{A}}^{(2)}}$ (see Definition~\ref{DCongOpC} for the general case).
Let us also recall, from Definition~\ref{DCongOpC}, that, if $\Phi\subseteq \mathrm{Pth}_{\boldsymbol{\mathcal{A}}^{(2)}}\times \mathrm{Pth}_{\boldsymbol{\mathcal{A}}^{(2)}}$, then  
$$
\mathrm{C}_{\mathbf{Pth}_{\boldsymbol{\mathcal{A}}^{(2)}}}\left(
\Phi
\right)
=\left(
\Phi\circ\Phi
\right)
\cup
\left(
\bigcup_{\gamma\in\Sigma^{\boldsymbol{\mathcal{A}}^{(2)}}_{\neq\lambda,s}}
\gamma^{\mathbf{Pth}_{\boldsymbol{\mathcal{A}}^{(2)}}}
\times
\gamma^{\mathbf{Pth}_{\boldsymbol{\mathcal{A}}^{(2)}}}
\left[
\Phi_{\mathrm{ar}(\gamma)}
\right]
\right)_{s\in S}.
$$
Moreover, for the family $(\mathrm{C}^{n}_{\mathbf{Pth}_{\boldsymbol{\mathcal{A}}^{(2)}}}(\Upsilon^{(1)}))_{n\in\mathbb{N}}$ in $\mathrm{Sub}(\mathrm{Pth}_{\boldsymbol{\mathcal{A}}^{(2)}}\times \mathrm{Pth}_{\boldsymbol{\mathcal{A}}^{(2)}})$, defined recursively, as follows:
\allowdisplaybreaks
\begin{align*}
\mathrm{C}^{0}_{\mathbf{Pth}_{\boldsymbol{\mathcal{A}}^{(2)}}}\left(
\Upsilon^{(1)}\right)
&=
\Upsilon^{(1)}\cup\left(\Upsilon^{(1)}\right)^{-1}\cup\Delta_{\mathrm{Pth}_{\boldsymbol{\mathcal{A}}^{(2)}}},\\
\mathrm{C}^{n+1}_{\mathbf{Pth}_{\boldsymbol{\mathcal{A}}^{(2)}}}\left(
\Upsilon^{(1)}\right)
&=\mathrm{C}_{\mathbf{Pth}_{\boldsymbol{\mathcal{A}}^{(2)}}}
\left(
\mathrm{C}^{n}_{\mathbf{Pth}_{\boldsymbol{\mathcal{A}}^{(2)}}}\left(
\Upsilon^{(1)}\right)\right),\, n\geq 0,
\end{align*}
we have, by Proposition~\ref{PCongOpC}, that 
$$\mathrm{C}^{\omega}_{\mathbf{Pth}_{\boldsymbol{\mathcal{A}}^{(2)}}}\left(
\Upsilon^{(1)}\right)
=\bigcup_{n\in\mathbb{N}}\mathrm{C}^{n}_{\mathbf{Pth}_{\boldsymbol{\mathcal{A}}^{(2)}}}\left(
\Upsilon^{(1)}
\right)
=\Upsilon^{[1]},$$ 
is equal to the smallest $\Sigma^{\boldsymbol{\mathcal{A}}^{(2)}}$-congruence on $\mathbf{Pth}_{\boldsymbol{\mathcal{A}}^{(2)}}$ containing $\Upsilon^{(1)}$.
\end{remark}

In the following lemmas we investigate the $\Sigma^{\boldsymbol{\mathcal{A}}^{(2)}}$-congruence $\Upsilon^{[1]}$ on $\mathbf{Pth}_{\boldsymbol{\mathcal{A}}^{(2)}}$. 

We start with two results that will be of great importance later on. They state that the results obtained in Corollaries~\ref{LDUps} and~\ref{CDUps}  for the relation $\Upsilon^{(1)}$ also hold for the $\Sigma^{\boldsymbol{\mathcal{A}}^{(2)}}$-congruence $\Upsilon^{[1]}$.

In the following lemma we prove that two second-order paths of the same sort that are $\Upsilon^{[1]}$-related must have the same length, the same $([1],2)$-source and the same $([1],2)$-target.

\begin{restatable}{lemma}{LDUpsCong}
\label{LDUpsCong} Let $s$ be a sort in $S$ and $(\mathfrak{P}^{(2)},\mathfrak{Q}^{(2)})\in \Upsilon^{[1]}_{s}$. The following statements hold.
\begin{enumerate}
\item[(i)] $\bb{\mathfrak{P}^{(2)}}=\bb{\mathfrak{Q}^{(2)}}$;
\item[(ii)] $\mathrm{sc}_{s}^{([1],2)}(\mathfrak{P}^{(2)})=\mathrm{sc}_{s}^{([1],2)}(\mathfrak{Q}^{(2)})$;
\item[(iii)] $\mathrm{tg}_{s}^{([1],2)}(\mathfrak{P}^{(2)})=\mathrm{tg}_{s}^{([1],2)}(\mathfrak{Q}^{(2)})$.
\end{enumerate}
\end{restatable}
\begin{proof} 
We recall from Remark~\ref{RDUpsCong} that 
$$\mathrm{C}^{\omega}_{\mathbf{Pth}_{\boldsymbol{\mathcal{A}}^{(2)}}}\left(
\Upsilon^{(1)}\right)
=\bigcup_{n\in\mathbb{N}}\mathrm{C}^{n}_{\mathbf{Pth}_{\boldsymbol{\mathcal{A}}^{(2)}}}\left(
\Upsilon^{(1)}
\right)
=\Upsilon^{[1]}.$$ 

We prove the statement by induction on  $n\in\mathbb{N}$.

\textsf{Base step of the induction.}

Let us recall from Remark~\ref{RDUpsCong} that 
\[
\mathrm{C}^{0}_{\mathbf{Pth}_{\boldsymbol{\mathcal{A}}^{(2)}}}\left(
\Upsilon^{(1)}\right)
=
\Upsilon^{(1)}\cup\left(\Upsilon^{(1)}\right)^{-1}\cup\Delta_{\mathrm{Pth}_{\boldsymbol{\mathcal{A}}^{(2)}}}.
\]

The statement trivially holds for pairs $(\mathfrak{P}^{(2)},\mathfrak{Q}^{(2)})\in \Delta_{\mathrm{Pth}_{\boldsymbol{\mathcal{A}}^{(2)},s}}$. 

If the pair $(\mathfrak{P}^{(2)},\mathfrak{Q}^{(2)})$ is in $\Upsilon^{(1)}_{s}$, the statement follows from Lemma~\ref{LDUps}. In case $(\mathfrak{P}^{(2)},\mathfrak{Q}^{(2)})$ is a pair in $(\Upsilon^{(1)}_{s})^{-1}$, we reason in a similar way.

This completes the base case.

\textsf{Inductive step of the induction.}

Assume that the statement holds for $n\in\mathbb{N}$, i.e., for every sort $s\in S$ and every pair of paths $(\mathfrak{P}^{(2)},\mathfrak{Q}^{(2)})$ in $\mathrm{C}^{n}_{\mathbf{Pth}_{\boldsymbol{\mathcal{A}}^{(2)}}}\left(
\Upsilon^{(1)}\right)$ it holds that 
\begin{enumerate}
\item[(i)] $\bb{\mathfrak{P}^{(2)}}=\bb{\mathfrak{Q}^{(2)}}$;
\item[(ii)] $\mathrm{sc}_{s}^{([1],2)}(\mathfrak{P}^{(2)})=\mathrm{sc}_{s}^{([1],2)}(\mathfrak{Q}^{(2)})$;
\item[(iii)] $\mathrm{tg}_{s}^{([1],2)}(\mathfrak{P}^{(2)})=\mathrm{tg}_{s}^{([1],2)}(\mathfrak{Q}^{(2)})$.
\end{enumerate}

We now prove the statement for $n+1$. Let $s$ be a sort in $S$ and let $(\mathfrak{P}^{(2)},\mathfrak{Q}^{(2)})$ be a pair of second-order paths in $\mathrm{Pth}_{\boldsymbol{\mathcal{A}}^{(2)},s}$ such that $(\mathfrak{P}^{(2)},\mathfrak{Q}^{(2)})\in \mathrm{C}^{n+1}_{\mathbf{Pth}_{\boldsymbol{\mathcal{A}}^{(2)}}}\left(
\Upsilon^{(1)}\right)$. Let us recall from Remark~\ref{RDUpsCong} that 
\allowdisplaybreaks
\begin{multline*}
\mathrm{C}^{n+1}_{\mathbf{Pth}_{\boldsymbol{\mathcal{A}}^{(2)}}}\left(
\Upsilon^{(1)}
\right)
=\left(
\mathrm{C}^{n}_{\mathbf{Pth}_{\boldsymbol{\mathcal{A}}^{(2)}}}\left(
\Upsilon^{(1)}
\right)
\circ
\mathrm{C}^{n}_{\mathbf{Pth}_{\boldsymbol{\mathcal{A}}^{(2)}}}\left(
\Upsilon^{(1)}
\right)
\right)
\cup
\\
\left(
\bigcup_{\gamma\in\Sigma^{\boldsymbol{\mathcal{A}}^{(2)}}_{\neq\lambda,s}}
\gamma^{\mathbf{Pth}_{\boldsymbol{\mathcal{A}}^{(2)}}}
\times
\gamma^{\mathbf{Pth}_{\boldsymbol{\mathcal{A}}^{(2)}}}
\left[
\mathrm{C}^{n}_{\mathbf{Pth}_{\boldsymbol{\mathcal{A}}^{(2)}}}\left(
\Upsilon^{(1)}
\right)_{\mathrm{ar}(\gamma)}
\right]
\right)_{s\in S}.
\end{multline*}

Then either (1),  $(\mathfrak{P}^{(2)},\mathfrak{Q}^{(2)})$ is a pair in $\mathrm{C}^{n}_{\mathbf{Pth}_{\boldsymbol{\mathcal{A}}^{(2)}}}(
\Upsilon^{(1)}
)
\circ
\mathrm{C}^{n}_{\mathbf{Pth}_{\boldsymbol{\mathcal{A}}^{(2)}}}(
\Upsilon^{(1)}
)$ or (2),  $(\mathfrak{P}^{(2)},\mathfrak{Q}^{(2)})$ is a pair in $\gamma^{\mathbf{Pth}_{\boldsymbol{\mathcal{A}}^{(2)}}}
\times
\gamma^{\mathbf{Pth}_{\boldsymbol{\mathcal{A}}^{(2)}}}
[
\mathrm{C}^{n}_{\mathbf{Pth}_{\boldsymbol{\mathcal{A}}^{(2)}}}(
\Upsilon^{(1)}
)_{\mathrm{ar}(\gamma)}]$ for some operation symbol $\gamma\in \Sigma^{\boldsymbol{\mathcal{A}}^{(2)}}_{\neq\lambda, s}$.

If (1), then there exists a second-order path $\mathfrak{R}^{(2)}\in\mathrm{Pth}_{\boldsymbol{\mathcal{A}}^{(2)},s}$ for which $(\mathfrak{P}^{(2)}, \mathfrak{R}^{(2)})$ and $(\mathfrak{R}^{(2)}, \mathfrak{Q}^{(2)})$ belong to $\mathrm{C}^{n}_{\mathbf{Pth}_{\boldsymbol{\mathcal{A}}^{(2)}}}(
\Upsilon^{(1)})$.  Since $(\mathfrak{P}^{(2)}, \mathfrak{R}^{(2)})$ is in $\mathrm{C}^{n}_{\mathbf{Pth}_{\boldsymbol{\mathcal{A}}^{(2)}}}(
\Upsilon^{(1)})_{s}$ then, by induction, we have that 
\begin{enumerate}
\item[(i)] $\bb{\mathfrak{P}^{(2)}}=\bb{\mathfrak{R}^{(2)}}$;
\item[(ii)] $\mathrm{sc}_{s}^{([1],2)}(\mathfrak{P}^{(2)})=\mathrm{sc}_{s}^{([1],2)}(\mathfrak{R}^{(2)})$;
\item[(iii)] $\mathrm{tg}_{s}^{([1],2)}(\mathfrak{P}^{(2)})=\mathrm{tg}_{s}^{([1],2)}(\mathfrak{R}^{(2)})$.
\end{enumerate}

Since $(\mathfrak{R}^{(2)}, \mathfrak{Q}^{(2)})$ is in $\mathrm{C}^{n}_{\mathbf{Pth}_{\boldsymbol{\mathcal{A}}^{(2)}}}(
\Upsilon^{(1)})_{s}$ then, by induction, we have that 
\begin{enumerate}
\item[(i)] $\bb{\mathfrak{R}^{(2)}}=\bb{\mathfrak{Q}^{(2)}}$;
\item[(ii)] $\mathrm{sc}_{s}^{([1],2)}(\mathfrak{R}^{(2)})=\mathrm{sc}_{s}^{([1],2)}(\mathfrak{Q}^{(2)})$;
\item[(iii)] $\mathrm{tg}_{s}^{([1],2)}(\mathfrak{R}^{(2)})=\mathrm{tg}_{s}^{([1],2)}(\mathfrak{Q}^{(2)})$.
\end{enumerate}

All in all, we conclude that 
\begin{enumerate}
\item[(i)] $\bb{\mathfrak{P}^{(2)}}=\bb{\mathfrak{Q}^{(2)}}$;
\item[(ii)] $\mathrm{sc}_{s}^{([1],2)}(\mathfrak{P}^{(2)})=\mathrm{sc}_{s}^{([1],2)}(\mathfrak{Q}^{(2)})$;
\item[(iii)] $\mathrm{tg}_{s}^{([1],2)}(\mathfrak{P}^{(2)})=\mathrm{tg}_{s}^{([1],2)}(\mathfrak{Q}^{(2)})$.
\end{enumerate}

This completes the Case~(1).

If~(2), then there exists a unique word $\mathbf{s}\in S^{\star}-\{\lambda\}$, a unique $s\in S$, a unique operation symbol $\gamma\in \Sigma^{\boldsymbol{\mathcal{A}}^{(2)}}_{\mathbf{s},s}$ and a unique family of pairs $((\mathfrak{P}^{(2)}_{j}, \mathfrak{Q}^{(2)}_{j}))_{j\in\bb{\mathbf{s}}}$ in $\mathrm{C}^{n}_{\mathbf{Pth}_{\boldsymbol{\mathcal{A}}^{(2)}}}(\Upsilon^{(1)})_{\mathbf{s}}$ in the domain of $\gamma$ for which 
\[
\left(
\mathfrak{P}^{(2)},
\mathfrak{Q}^{(2)}
\right)
=
\left(
\gamma^{\mathbf{Pth}_{\boldsymbol{\mathcal{A}}^{(2)}}}\left(\left(\mathfrak{P}^{(2)}_{j}\right)_{j\in\bb{\mathbf{s}}}\right),
\gamma^{\mathbf{Pth}_{\boldsymbol{\mathcal{A}}^{(2)}}}\left(\left(\mathfrak{Q}^{(2)}_{j}\right)_{j\in\bb{\mathbf{s}}}\right)
\right).
\]

We will distinguish the following cases according to the different possibilities for the operation symbol $\gamma\in \Sigma^{\boldsymbol{\mathcal{A}}^{(2)}}_{\mathbf{s},s}$. Note that either (2.1) $\gamma$ is an operation symbol $\sigma\in \Sigma_{\mathbf{s},s}$; or (2.2) $\gamma$ is the operation symbol of $0$-source $\mathrm{sc}^{0}_{s}\in \Sigma^{\boldsymbol{\mathcal{A}}}_{s,s}$; or (2.3) $\gamma$ is the operation symbol of $0$-target $\mathrm{tg}^{0}_{s}\in \Sigma^{\boldsymbol{\mathcal{A}}}_{s,s}$; or (2.4) $\gamma$ is the operation symbol of $0$-composition $\circ^{0}_{s}\in \Sigma^{\boldsymbol{\mathcal{A}}}_{ss,s}$; or (2.5) $\gamma$ is the operation symbol of $1$-source $\mathrm{sc}^{1}_{s}\in \Sigma^{\boldsymbol{\mathcal{A}}^{(2)}}_{s,s}$; or (2.6) $\gamma$ is the operation symbol of $1$-target $\mathrm{tg}^{1}_{s}\in \Sigma^{\boldsymbol{\mathcal{A}}^{(2)}}_{s,s}$; or (2.7) $\gamma$ is the operation symbol of $1$-composition $\circ^{1}_{s}\in \Sigma^{\boldsymbol{\mathcal{A}}^{(2)}}_{ss,s}$.

\textsf{Case~(2.1)} $\gamma$ is an operation symbol $\sigma \in \Sigma_{\mathbf{s},s}$. 

Let $\mathbf{s}$ be a word in $S^{\star}-\{\lambda\}$, let $\sigma$ be an operation symbol in $\Sigma_{\mathbf{s},s}$ and let $((\mathfrak{P}^{(2)}_{j}, \mathfrak{Q}^{(2)}_{j}))_{j\in\bb{\mathbf{s}}}$ be the family of pairs in $\mathrm{C}^{n}_{\mathbf{Pth}_{\boldsymbol{\mathcal{A}}^{(2)}}}(\Upsilon^{(1)})_{\mathbf{s}}$ for which 
\[
\left(
\mathfrak{P}^{(2)},
\mathfrak{Q}^{(2)}
\right)
=
\left(
\sigma^{\mathbf{Pth}_{\boldsymbol{\mathcal{A}}^{(2)}}}\left(\left(\mathfrak{P}^{(2)}_{j}\right)_{j\in\bb{\mathbf{s}}}\right),
\sigma^{\mathbf{Pth}_{\boldsymbol{\mathcal{A}}^{(2)}}}\left(\left(\mathfrak{Q}^{(2)}_{j}\right)_{j\in\bb{\mathbf{s}}}\right)
\right).
\]

Let us note that, for every $j\in\bb{\mathbf{s}}$, since $(\mathfrak{P}^{(2)}_{j}, \mathfrak{Q}^{(2)}_{j})$ is a pair in $\mathrm{C}^{n}_{\mathbf{Pth}_{\boldsymbol{\mathcal{A}}^{(2)}}}(\Upsilon^{(1)})_{s_{j}}$, we have, by induction, that 
\begin{enumerate}
\item[(i)] $\bb{\mathfrak{P}^{(2)}_{j}}=\bb{\mathfrak{Q}^{(2)}_{j}}$;
\item[(ii)] $\mathrm{sc}_{s_{j}}^{([1],2)}(\mathfrak{P}^{(2)}_{j})=\mathrm{sc}_{s_{j}}^{([1],2)}(\mathfrak{Q}^{(2)}_{j})$;
\item[(iii)] $\mathrm{tg}_{s_{j}}^{([1],2)}(\mathfrak{P}^{(2)}_{j})=\mathrm{tg}_{s_{j}}^{([1],2)}(\mathfrak{Q}^{(2)}_{j})$.
\end{enumerate}

Note that the following chain of equalities hold
\allowdisplaybreaks
\begin{align*}
\left\lvert
\sigma^{\mathbf{Pth}_{\boldsymbol{\mathcal{A}}^{(2)}}}\left(\left(\mathfrak{P}^{(2)}_{j}\right)_{j\in\bb{\mathbf{s}}}\right)
\right\rvert
&=
\sum_{j\in\bb{\mathbf{s}}}\left\lvert
\mathfrak{P}^{(2)}_{j}
\right\rvert
\tag{1}
\\&=
\sum_{j\in\bb{\mathbf{s}}}\left\lvert
\mathfrak{Q}^{(2)}_{j}
\right\rvert
\tag{2}
\\&=
\left\lvert
\sigma^{\mathbf{Pth}_{\boldsymbol{\mathcal{A}}^{(2)}}}\left(\left(\mathfrak{Q}^{(2)}_{j}\right)_{j\in\bb{\mathbf{s}}}\right)
\right\rvert.
\tag{3}
\end{align*}

In the just stated chain of equalities, the first equality follows from Claim~\ref{CDPthSigma}; the second equality follows by induction; finally, the third equality follows from Claim~\ref{CDPthSigma}.

Regarding the $([1],2)$-source, the following chain of equalities hold
\allowdisplaybreaks
\begin{align*}
\mathrm{sc}^{([1],2)}_{s}\left(
\sigma^{\mathbf{Pth}_{\boldsymbol{\mathcal{A}}^{(2)}}}\left(\left(\mathfrak{P}^{(2)}_{j}\right)_{j\in\bb{\mathbf{s}}}\right)
\right)
&=
\sigma^{[\mathbf{PT}_{\boldsymbol{\mathcal{A}}}]}\left(
\left(
\mathrm{sc}^{([1],2)}_{s_{j}}\left(
\mathfrak{P}^{(2)}_{j}
\right)
\right)_{j\in\bb{\mathbf{s}}}
\right)
\tag{1}
\\&=
\sigma^{[\mathbf{PT}_{\boldsymbol{\mathcal{A}}}]}\left(
\left(
\mathrm{sc}^{([1],2)}_{s_{j}}\left(
\mathfrak{Q}^{(2)}_{j}
\right)
\right)_{j\in\bb{\mathbf{s}}}
\right)
\tag{2}
\\&=
\mathrm{sc}^{([1],2)}_{s}\left(
\sigma^{\mathbf{Pth}_{\boldsymbol{\mathcal{A}}^{(2)}}}\left(\left(\mathfrak{Q}^{(2)}_{j}\right)_{j\in\bb{\mathbf{s}}}\right)
\right).
\tag{3}
\end{align*}

In the just stated chain of equalities, the first equality follows from Claim~\ref{CDPthSigma}; the second equality follows by induction; finally, the third equality follows from Claim~\ref{CDPthSigma}.

Regarding the $([1],2)$-target, the following chain of equalities hold
\allowdisplaybreaks
\begin{align*}
\mathrm{tg}^{([1],2)}_{s}\left(
\sigma^{\mathbf{Pth}_{\boldsymbol{\mathcal{A}}^{(2)}}}\left(\left(\mathfrak{P}^{(2)}_{j}\right)_{j\in\bb{\mathbf{s}}}\right)
\right)
&=
\sigma^{[\mathbf{PT}_{\boldsymbol{\mathcal{A}}}]}\left(
\left(
\mathrm{tg}^{([1],2)}_{s_{j}}\left(
\mathfrak{P}^{(2)}_{j}
\right)
\right)_{j\in\bb{\mathbf{s}}}
\right)
\tag{1}
\\&=
\sigma^{[\mathbf{PT}_{\boldsymbol{\mathcal{A}}}]}\left(
\left(
\mathrm{tg}^{([1],2)}_{s_{j}}\left(
\mathfrak{Q}^{(2)}_{j}
\right)
\right)_{j\in\bb{\mathbf{s}}}
\right)
\tag{2}
\\&=
\mathrm{tg}^{([1],2)}_{s}\left(
\sigma^{\mathbf{Pth}_{\boldsymbol{\mathcal{A}}^{(2)}}}\left(\left(\mathfrak{Q}^{(2)}_{j}\right)_{j\in\bb{\mathbf{s}}}\right)
\right).
\tag{3}
\end{align*}

In the just stated chain of equalities, the first equality follows from Claim~\ref{CDPthSigma}; the second equality follows by induction; finally, the third equality follows from Claim~\ref{CDPthSigma}.

This completes the Case~(2.1).

\textsf{Case~(2.2)} $\gamma$ is the operation symbol $\mathrm{sc}^{0}_{s} \in \Sigma^{\boldsymbol{\mathcal{A}}}_{s,s}$. 

Let $(\mathfrak{P}'^{(2)}, \mathfrak{Q}'^{(2)})$ be the pairs in $\mathrm{C}^{n}_{\mathbf{Pth}_{\boldsymbol{\mathcal{A}}^{(2)}}}(\Upsilon^{(1)})_{s}$ for which 
\[
\left(
\mathfrak{P}^{(2)},
\mathfrak{Q}^{(2)}
\right)
=
\left(
\mathrm{sc}^{0\mathbf{Pth}_{\boldsymbol{\mathcal{A}}^{(2)}}}_{s}\left(\mathfrak{P}'^{(2)}\right),
\mathrm{sc}^{0\mathbf{Pth}_{\boldsymbol{\mathcal{A}}^{(2)}}}_{s}\left(\mathfrak{Q}'^{(2)}\right)
\right).
\]

Let us note that, since $(\mathfrak{P}'^{(2)}, \mathfrak{Q}'^{(2)})$ is a pair in $\mathrm{C}^{n}_{\mathbf{Pth}_{\boldsymbol{\mathcal{A}}^{(2)}}}(\Upsilon^{(1)})_{s}$, we have, by induction, that 
\begin{enumerate}
\item[(i)] $\bb{\mathfrak{P}'^{(2)}_{j}}=\bb{\mathfrak{Q}'^{(2)}_{j}}$;
\item[(ii)] $\mathrm{sc}_{s}^{([1],2)}(\mathfrak{P}'^{(2)})=\mathrm{sc}_{s}^{([1],2)}(\mathfrak{Q}'^{(2)})$;
\item[(iii)] $\mathrm{tg}_{s}^{([1],2)}(\mathfrak{P}'^{(2)})=\mathrm{tg}_{s}^{([1],2)}(\mathfrak{Q}'^{(2)})$.
\end{enumerate}

Note that in virtue of Claim~\ref{CDPthCatAlgScZ}, we have that 
\[
\left\lvert
\mathrm{sc}^{0\mathbf{Pth}_{\boldsymbol{\mathcal{A}}^{(2)}}}_{s}\left(\mathfrak{P}'^{(2)}\right)
\right\rvert
=
0
=
\left\lvert
\mathrm{sc}^{0\mathbf{Pth}_{\boldsymbol{\mathcal{A}}^{(2)}}}_{s}\left(\mathfrak{Q}'^{(2)}\right)
\right\rvert.
\]

Regarding the $([1],2)$-source, the following chain of equalities hold
\allowdisplaybreaks
\begin{align*}
\mathrm{sc}^{([1],2)}_{s}\left(
\mathrm{sc}^{0\mathbf{Pth}_{\boldsymbol{\mathcal{A}}^{(2)}}}_{s}\left(\mathfrak{P}'^{(2)}\right)
\right)
&=
\mathrm{sc}^{0[\mathbf{PT}_{\boldsymbol{\mathcal{A}}}]}_{s}\left(
\mathrm{sc}_{s}^{([1],2)}\left(\mathfrak{P}'^{(2)}\right)
\right)
\tag{1}
\\&=
\mathrm{sc}^{0[\mathbf{PT}_{\boldsymbol{\mathcal{A}}}]}_{s}\left(
\mathrm{sc}_{s}^{([1],2)}\left(\mathfrak{Q}'^{(2)}\right)
\right)
\tag{2}
\\&=
\mathrm{sc}^{([1],2)}_{s}\left(
\mathrm{sc}^{0\mathbf{Pth}_{\boldsymbol{\mathcal{A}}^{(2)}}}_{s}\left(\mathfrak{Q}'^{(2)}\right)
\right).
\tag{3}
\end{align*} 

In the just stated chain of equalities, the first equality follows from Proposition~\ref{PDUCatHom}; the second equality follows by induction; finally, the last equality follows from Proposition~\ref{PDUCatHom}.

Regarding the $([1],2)$-target, the following chain of equalities hold
\allowdisplaybreaks
\begin{align*}
\mathrm{tg}^{([1],2)}_{s}\left(
\mathrm{sc}^{0\mathbf{Pth}_{\boldsymbol{\mathcal{A}}^{(2)}}}_{s}\left(\mathfrak{P}'^{(2)}\right)
\right)
&=
\mathrm{tg}^{0[\mathbf{PT}_{\boldsymbol{\mathcal{A}}}]}_{s}\left(
\mathrm{sc}_{s}^{([1],2)}\left(\mathfrak{P}'^{(2)}\right)
\right)
\tag{1}
\\&=
\mathrm{tg}^{0[\mathbf{PT}_{\boldsymbol{\mathcal{A}}}]}_{s}\left(
\mathrm{sc}_{s}^{([1],2)}\left(\mathfrak{Q}'^{(2)}\right)
\right)
\tag{2}
\\&=
\mathrm{tg}^{([1],2)}_{s}\left(
\mathrm{sc}^{0\mathbf{Pth}_{\boldsymbol{\mathcal{A}}^{(2)}}}_{s}\left(\mathfrak{Q}'^{(2)}\right)
\right).
\tag{3}
\end{align*} 

In the just stated chain of equalities, the first equality follows from Proposition~\ref{PDUCatHom}; the second equality follows by induction; finally, the last equality follows from Proposition~\ref{PDUCatHom}.

This completes the Case~(2.2).

\textsf{Case~(2.3)} $\gamma$ is the operation symbol $\mathrm{tg}^{0}_{s} \in \Sigma^{\boldsymbol{\mathcal{A}}}_{s,s}$. 

Let $(\mathfrak{P}'^{(2)}, \mathfrak{Q}'^{(2)})$ be the pairs in $\mathrm{C}^{n}_{\mathbf{Pth}_{\boldsymbol{\mathcal{A}}^{(2)}}}(\Upsilon^{(1)})_{s}$ for which 
\[
\left(
\mathfrak{P}^{(2)},
\mathfrak{Q}^{(2)}
\right)
=
\left(
\mathrm{tg}^{0\mathbf{Pth}_{\boldsymbol{\mathcal{A}}^{(2)}}}_{s}\left(\mathfrak{P}'^{(2)}\right),
\mathrm{tg}^{0\mathbf{Pth}_{\boldsymbol{\mathcal{A}}^{(2)}}}_{s}\left(\mathfrak{Q}'^{(2)}\right)
\right).
\]

This case can be proven by a similar reasoning than that used for Case~(2.2).

This completes the Case~(2.3).

\textsf{Case~(2.4)} $\gamma$ is the operation symbol $\circ^{0}_{s} \in \Sigma^{\boldsymbol{\mathcal{A}}}_{ss,s}$. 

Let $(\mathfrak{P}'^{(2)}, \mathfrak{Q}'^{(2)})$ and $(\mathfrak{P}''^{(2)}, \mathfrak{Q}''^{(2)})$ be the pairs in $\mathrm{C}^{n}_{\mathbf{Pth}_{\boldsymbol{\mathcal{A}}^{(2)}}}(\Upsilon^{(1)})_{s}$ satisfying 
\begin{align*}
\mathrm{sc}^{(0,2)}_{s}\left(
\mathfrak{P}''^{(2)}
\right)
&=
\mathrm{tg}^{(0,2)}_{s}\left(
\mathfrak{P}'^{(2)}
\right);
&
\mathrm{sc}^{(0,2)}_{s}\left(
\mathfrak{Q}''^{(2)}
\right)
&=
\mathrm{tg}^{(0,2)}_{s}\left(
\mathfrak{Q}'^{(2)}
\right),
\end{align*}
for which 
\[
\left(
\mathfrak{P}^{(2)},
\mathfrak{Q}^{(2)}
\right)
=
\left(
\mathfrak{P}''^{(2)}\circ^{0\mathbf{Pth}_{\boldsymbol{\mathcal{A}}^{(2)}}}_{s} \mathfrak{P}'^{(2)},
\mathfrak{Q}''^{(2)}\circ^{0\mathbf{Pth}_{\boldsymbol{\mathcal{A}}^{(2)}}}_{s} \mathfrak{Q}'^{(2)}
\right).
\]

Let us note that, since $(\mathfrak{P}'^{(2)}, \mathfrak{Q}'^{(2)})$ is a pair in $\mathrm{C}^{n}_{\mathbf{Pth}_{\boldsymbol{\mathcal{A}}^{(2)}}}(\Upsilon^{(1)})_{s}$, we have, by induction, that 
\begin{enumerate}
\item[(i)] $\bb{\mathfrak{P}'^{(2)}}=\bb{\mathfrak{Q}'^{(2)}}$;
\item[(ii)] $\mathrm{sc}_{s}^{([1],2)}(\mathfrak{P}'^{(2)})=\mathrm{sc}_{s}^{([1],2)}(\mathfrak{Q}'^{(2)})$;
\item[(iii)] $\mathrm{tg}_{s}^{([1],2)}(\mathfrak{P}'^{(2)})=\mathrm{tg}_{s}^{([1],2)}(\mathfrak{Q}'^{(2)})$.
\end{enumerate}

Since $(\mathfrak{P}''^{(2)}, \mathfrak{Q}''^{(2)})$ is a pair in $\mathrm{C}^{n}_{\mathbf{Pth}_{\boldsymbol{\mathcal{A}}^{(2)}}}(\Upsilon^{(1)})_{s}$, we have, by induction, that 
\begin{enumerate}
\item[(i)] $\bb{\mathfrak{P}''^{(2)}}=\bb{\mathfrak{Q}''^{(2)}}$;
\item[(ii)] $\mathrm{sc}_{s}^{([1],2)}(\mathfrak{P}''^{(2)})=\mathrm{sc}_{s}^{([1],2)}(\mathfrak{Q}''^{(2)})$;
\item[(iii)] $\mathrm{tg}_{s}^{([1],2)}(\mathfrak{P}''^{(2)})=\mathrm{tg}_{s}^{([1],2)}(\mathfrak{Q}''^{(2)})$.
\end{enumerate}

Note that the following chain of equalities holds
\allowdisplaybreaks
\begin{align*}
\left\lvert
\mathfrak{P}''^{(2)}\circ^{0\mathbf{Pth}_{\boldsymbol{\mathcal{A}}^{(2)}}}_{s} \mathfrak{P}'^{(2)}
\right\rvert
&=
\left\lvert
\mathfrak{P}''^{(2)}
\right\rvert
+
\left\lvert
\mathfrak{P}'^{(2)}
\right\rvert
\tag{1}
\\&=
\left\lvert
\mathfrak{Q}''^{(2)}
\right\rvert
+
\left\lvert
\mathfrak{Q}'^{(2)}
\right\rvert
\tag{2}
\\&=
\left\lvert
\mathfrak{Q}''^{(2)}\circ^{0\mathbf{Pth}_{\boldsymbol{\mathcal{A}}^{(2)}}}_{s} \mathfrak{Q}'^{(2)}
\right\rvert.
\tag{3}
\end{align*}

In the just stated chain of equalities, the first equality follows from Claim~\ref{CDPthCatAlgCompZ}; the second equality follows by induction; finally, the last equality follows from Claim~\ref{CDPthCatAlgCompZ}.

Regarding the $([1],2)$-source, the following chain of equalities hold
\begin{flushleft}
$\mathrm{sc}^{([1],2)}_{s}\left(
\mathfrak{P}''^{(2)}\circ^{0\mathbf{Pth}_{\boldsymbol{\mathcal{A}}^{(2)}}}_{s} \mathfrak{P}'^{(2)}
\right)$
\allowdisplaybreaks
\begin{align*}
\qquad
&=
\left(
\mathrm{sc}_{s}^{([1],2)}\left(
\mathfrak{P}''^{(2)}
\right)
\right)
\circ^{0[\mathbf{PT}_{\boldsymbol{\mathcal{A}}}]}_{s}
\left(
\mathrm{sc}_{s}^{([1],2)}\left(
\mathfrak{P}'^{(2)}
\right)
\right)
\tag{1}
\\&=
\left(
\mathrm{sc}_{s}^{([1],2)}\left(
\mathfrak{Q}''^{(2)}
\right)
\right)
\circ^{0[\mathbf{PT}_{\boldsymbol{\mathcal{A}}}]}_{s}
\left(
\mathrm{sc}_{s}^{([1],2)}\left(
\mathfrak{Q}'^{(2)}
\right)
\right)
\tag{2}
\\&=
\mathrm{sc}^{([1],2)}_{s}\left(
\mathfrak{Q}''^{(2)}\circ^{0\mathbf{Pth}_{\boldsymbol{\mathcal{A}}^{(2)}}}_{s} \mathfrak{Q}'^{(2)}
\right).
\tag{3}
\end{align*} 
\end{flushleft}

In the just stated chain of equalities, the first equality follows from Claim~\ref{CDPthCatAlgCompZ}; the second equality follows by induction; finally, the last equality follows from Claim~\ref{CDPthCatAlgCompZ}.

Regarding the $([1],2)$-target, the following chain of equalities hold
\begin{flushleft}
$\mathrm{tg}^{([1],2)}_{s}\left(
\mathfrak{P}''^{(2)}\circ^{0\mathbf{Pth}_{\boldsymbol{\mathcal{A}}^{(2)}}}_{s} \mathfrak{P}'^{(2)}
\right)$
\allowdisplaybreaks
\begin{align*}
\qquad
&=
\left(
\mathrm{tg}_{s}^{([1],2)}\left(
\mathfrak{P}''^{(2)}
\right)
\right)
\circ^{0[\mathbf{PT}_{\boldsymbol{\mathcal{A}}}]}_{s}
\left(
\mathrm{tg}_{s}^{([1],2)}\left(
\mathfrak{P}'^{(2)}
\right)
\right)
\tag{1}
\\&=
\left(
\mathrm{tg}_{s}^{([1],2)}\left(
\mathfrak{Q}''^{(2)}
\right)
\right)
\circ^{0[\mathbf{PT}_{\boldsymbol{\mathcal{A}}}]}_{s}
\left(
\mathrm{tg}_{s}^{([1],2)}\left(
\mathfrak{Q}'^{(2)}
\right)
\right)
\tag{2}
\\&=
\mathrm{tg}^{([1],2)}_{s}\left(
\mathfrak{Q}''^{(2)}\circ^{0\mathbf{Pth}_{\boldsymbol{\mathcal{A}}^{(2)}}}_{s} \mathfrak{Q}'^{(2)}
\right).
\tag{3}
\end{align*} 
\end{flushleft}

In the just stated chain of equalities, the first equality follows from Claim~\ref{CDPthCatAlgCompZ}; the second equality follows by induction; finally, the last equality follows from Claim~\ref{CDPthCatAlgCompZ}.

This completes the Case~(2.4).

\textsf{Case~(2.5)} $\gamma$ is the operation symbol $\mathrm{sc}^{1}_{s} \in \Sigma^{\boldsymbol{\mathcal{A}}}_{s,s}$. 

Let $(\mathfrak{P}'^{(2)}, \mathfrak{Q}'^{(2)})$ be the pairs in $\mathrm{C}^{n}_{\mathbf{Pth}_{\boldsymbol{\mathcal{A}}^{(2)}}}(\Upsilon^{(1)})_{s}$ for which 
\[
\left(
\mathfrak{P}^{(2)},
\mathfrak{Q}^{(2)}
\right)
=
\left(
\mathrm{sc}^{1\mathbf{Pth}_{\boldsymbol{\mathcal{A}}^{(2)}}}_{s}\left(\mathfrak{P}'^{(2)}\right),
\mathrm{sc}^{1\mathbf{Pth}_{\boldsymbol{\mathcal{A}}^{(2)}}}_{s}\left(\mathfrak{Q}'^{(2)}\right)
\right).
\]

Let us note that, since $(\mathfrak{P}'^{(2)}, \mathfrak{Q}'^{(2)})$ is a pair in $\mathrm{C}^{n}_{\mathbf{Pth}_{\boldsymbol{\mathcal{A}}^{(2)}}}(\Upsilon^{(1)})_{s}$, we have, by induction, that 
\begin{enumerate}
\item[(i)] $\bb{\mathfrak{P}'^{(2)}_{j}}=\bb{\mathfrak{Q}'^{(2)}_{j}}$;
\item[(ii)] $\mathrm{sc}_{s}^{([1],2)}(\mathfrak{P}'^{(2)})=\mathrm{sc}_{s}^{([1],2)}(\mathfrak{Q}'^{(2)})$;
\item[(iii)] $\mathrm{tg}_{s}^{([1],2)}(\mathfrak{P}'^{(2)})=\mathrm{tg}_{s}^{([1],2)}(\mathfrak{Q}'^{(2)})$.
\end{enumerate}

Note that in virtue of Proposition~\ref{PDPthDCatAlg}, we have that 
\[
\left\lvert
\mathrm{sc}^{1\mathbf{Pth}_{\boldsymbol{\mathcal{A}}^{(2)}}}_{s}\left(\mathfrak{P}'^{(2)}\right)
\right\rvert
=
0
=
\left\lvert
\mathrm{sc}^{1\mathbf{Pth}_{\boldsymbol{\mathcal{A}}^{(2)}}}_{s}\left(\mathfrak{Q}'^{(2)}\right)
\right\rvert.
\]

Regarding the $([1],2)$-source, the following chain of equalities hold
\allowdisplaybreaks
\begin{align*}
\mathrm{sc}^{([1],2)}_{s}\left(
\mathrm{sc}^{1\mathbf{Pth}_{\boldsymbol{\mathcal{A}}^{(2)}}}_{s}\left(\mathfrak{P}'^{(2)}\right)
\right)
&=
\mathrm{sc}^{([1],2)}_{s}\left(
\mathrm{ip}^{(2,[1])\sharp}_{s}\left(
\mathrm{sc}^{([1],2)}_{s}\left(
\mathfrak{P}^{(2)}
\right)
\right)
\right)
\tag{1}
\\&=
\mathrm{sc}^{([1],2)}_{s}\left(
\mathfrak{P}^{(2)}
\right)
\tag{2}
\\&=
\mathrm{sc}^{([1],2)}_{s}\left(
\mathfrak{Q}^{(2)}
\right)
\tag{3}
\\&=
\mathrm{sc}^{([1],2)}_{s}\left(
\mathrm{ip}^{(2,[1])\sharp}_{s}\left(
\mathrm{sc}^{([1],2)}_{s}\left(
\mathfrak{Q}^{(2)}
\right)
\right)
\right)
\tag{4}
\\ &=
\mathrm{sc}^{([1],2)}_{s}\left(
\mathrm{sc}^{1\mathbf{Pth}_{\boldsymbol{\mathcal{A}}^{(2)}}}_{s}\left(\mathfrak{P}'^{(2)}\right)
\right).
\tag{5}
\end{align*} 

In the just stated chain of equalities, the first equality follows from Proposition~~\ref{PDPthDCatAlg}; the second equality follows from Proposition~\ref{PDBasicEq}; the third equality follows by induction; the fourth equality follows from Proposition~\ref{PDBasicEq}; finally, the last equality follows from Proposition~\ref{PDPthDCatAlg}.

Regarding the $([1],2)$-target, the following chain of equalities hold
\allowdisplaybreaks
\begin{align*}
\mathrm{tg}^{([1],2)}_{s}\left(
\mathrm{sc}^{1\mathbf{Pth}_{\boldsymbol{\mathcal{A}}^{(2)}}}_{s}\left(\mathfrak{P}'^{(2)}\right)
\right)
&=
\mathrm{tg}^{([1],2)}_{s}\left(
\mathrm{ip}^{(2,[1])\sharp}_{s}\left(
\mathrm{sc}^{([1],2)}_{s}\left(
\mathfrak{P}^{(2)}
\right)
\right)
\right)
\tag{1}
\\&=
\mathrm{sc}^{([1],2)}_{s}\left(
\mathfrak{P}^{(2)}
\right)
\tag{2}
\\&=
\mathrm{sc}^{([1],2)}_{s}\left(
\mathfrak{Q}^{(2)}
\right)
\tag{3}
\\&=
\mathrm{tg}^{([1],2)}_{s}\left(
\mathrm{ip}^{(2,[1])\sharp}_{s}\left(
\mathrm{sc}^{([1],2)}_{s}\left(
\mathfrak{Q}^{(2)}
\right)
\right)
\right)
\tag{4}
\\ &=
\mathrm{tg}^{([1],2)}_{s}\left(
\mathrm{sc}^{1\mathbf{Pth}_{\boldsymbol{\mathcal{A}}^{(2)}}}_{s}\left(\mathfrak{P}'^{(2)}\right)
\right).
\tag{5}
\end{align*} 

In the just stated chain of equalities, the first equality follows from Proposition~~\ref{PDPthDCatAlg}; the second equality follows from Proposition~\ref{PDBasicEq}; the third equality follows by induction; the fourth equality follows from Proposition~\ref{PDBasicEq}; finally, the last equality follows from Proposition~\ref{PDPthDCatAlg}.

This completes the Case~(2.5).

\textsf{Case~(2.6)} $\gamma$ is the operation symbol $\mathrm{tg}^{1}_{s} \in \Sigma^{\boldsymbol{\mathcal{A}}}_{s,s}$. 

Let $(\mathfrak{P}'^{(2)}, \mathfrak{Q}'^{(2)})$ be the pairs in $\mathrm{C}^{n}_{\mathbf{Pth}_{\boldsymbol{\mathcal{A}}^{(2)}}}(\Upsilon^{(1)})_{s}$ for which 
\[
\left(
\mathfrak{P}^{(2)},
\mathfrak{Q}^{(2)}
\right)
=
\left(
\mathrm{tg}^{1\mathbf{Pth}_{\boldsymbol{\mathcal{A}}^{(2)}}}_{s}\left(\mathfrak{P}'^{(2)}\right),
\mathrm{tg}^{1\mathbf{Pth}_{\boldsymbol{\mathcal{A}}^{(2)}}}_{s}\left(\mathfrak{Q}'^{(2)}\right)
\right).
\]

This case can be proven by a similar reasoning than that used for Case~(2.5).

This completes the Case~(2.6).

\textsf{Case~(2.7)} $\gamma$ is the operation symbol $\circ^{1}_{s} \in \Sigma^{\boldsymbol{\mathcal{A}}}_{ss,s}$. 

Let $(\mathfrak{P}'^{(2)}, \mathfrak{Q}'^{(2)})$ and $(\mathfrak{P}''^{(2)}, \mathfrak{Q}''^{(2)})$ be the pairs in $\mathrm{C}^{n}_{\mathbf{Pth}_{\boldsymbol{\mathcal{A}}^{(2)}}}(\Upsilon^{(1)})_{s}$ satisfying 
\begin{align*}
\mathrm{sc}^{(0,2)}_{s}\left(
\mathfrak{P}''^{(2)}
\right)
&=
\mathrm{tg}^{(0,2)}_{s}\left(
\mathfrak{P}'^{(2)}
\right);
&
\mathrm{sc}^{(0,2)}_{s}\left(
\mathfrak{Q}''^{(2)}
\right)
&=
\mathrm{tg}^{(0,2)}_{s}\left(
\mathfrak{Q}'^{(2)}
\right),
\end{align*}
for which 
\[
\left(
\mathfrak{P}^{(2)},
\mathfrak{Q}^{(2)}
\right)
=
\left(
\mathfrak{P}''^{(2)}\circ^{1\mathbf{Pth}_{\boldsymbol{\mathcal{A}}^{(2)}}}_{s} \mathfrak{P}'^{(2)},
\mathfrak{Q}''^{(2)}\circ^{1\mathbf{Pth}_{\boldsymbol{\mathcal{A}}^{(2)}}}_{s} \mathfrak{Q}'^{(2)}
\right).
\]

Let us note that, since $(\mathfrak{P}'^{(2)}, \mathfrak{Q}'^{(2)})$ is a pair in $\mathrm{C}^{n}_{\mathbf{Pth}_{\boldsymbol{\mathcal{A}}^{(2)}}}(\Upsilon^{(1)})_{s}$, we have, by induction, that 
\begin{enumerate}
\item[(i)] $\bb{\mathfrak{P}'^{(2)}}=\bb{\mathfrak{Q}'^{(2)}}$;
\item[(ii)] $\mathrm{sc}_{s}^{([1],2)}(\mathfrak{P}'^{(2)})=\mathrm{sc}_{s}^{([1],2)}(\mathfrak{Q}'^{(2)})$;
\item[(iii)] $\mathrm{tg}_{s}^{([1],2)}(\mathfrak{P}'^{(2)})=\mathrm{tg}_{s}^{([1],2)}(\mathfrak{Q}'^{(2)})$.
\end{enumerate}

Since $(\mathfrak{P}''^{(2)}, \mathfrak{Q}''^{(2)})$ is a pair in $\mathrm{C}^{n}_{\mathbf{Pth}_{\boldsymbol{\mathcal{A}}^{(2)}}}(\Upsilon^{(1)})_{s}$, we have, by induction, that 
\begin{enumerate}
\item[(i)] $\bb{\mathfrak{P}''^{(2)}}=\bb{\mathfrak{Q}''^{(2)}}$;
\item[(ii)] $\mathrm{sc}_{s}^{([1],2)}(\mathfrak{P}''^{(2)})=\mathrm{sc}_{s}^{([1],2)}(\mathfrak{Q}''^{(2)})$;
\item[(iii)] $\mathrm{tg}_{s}^{([1],2)}(\mathfrak{P}''^{(2)})=\mathrm{tg}_{s}^{([1],2)}(\mathfrak{Q}''^{(2)})$.
\end{enumerate}

Note that the following chain of equalities hold
\allowdisplaybreaks
\begin{align*}
\left\lvert
\mathfrak{P}''^{(2)}\circ^{1\mathbf{Pth}_{\boldsymbol{\mathcal{A}}^{(2)}}}_{s} \mathfrak{P}'^{(2)}
\right\rvert
&=
\left\lvert
\mathfrak{P}''^{(2)}
\right\rvert
+
\left\lvert
\mathfrak{P}'^{(2)}
\right\rvert
\tag{1}
\\&=
\left\lvert
\mathfrak{Q}''^{(2)}
\right\rvert
+
\left\lvert
\mathfrak{Q}'^{(2)}
\right\rvert
\tag{2}
\\&=
\left\lvert
\mathfrak{Q}''^{(2)}\circ^{1\mathbf{Pth}_{\boldsymbol{\mathcal{A}}^{(2)}}}_{s} \mathfrak{Q}'^{(2)}
\right\rvert.
\tag{3}
\end{align*}

In the just stated chain of equalities, the first equality follows from Proposition~\ref{PDPthDCatAlg}; the second equality follows by induction; finally, the last equality follows from Proposition~\ref{PDPthDCatAlg}.

Regarding the $([1],2)$-source, the following chain of equalities hold
\allowdisplaybreaks
\begin{align*}
\mathrm{sc}^{([1],2)}_{s}\left(
\mathfrak{P}''^{(2)}\circ^{1\mathbf{Pth}_{\boldsymbol{\mathcal{A}}^{(2)}}}_{s} \mathfrak{P}'^{(2)}
\right)
&=
\mathrm{sc}_{s}^{([1],2)}\left(
\mathfrak{P}'^{(2)}
\right)
\tag{1}
\\&=
\mathrm{sc}_{s}^{([1],2)}\left(
\mathfrak{Q}'^{(2)}
\right)
\tag{2}
\\&=
\mathrm{sc}^{([1],2)}_{s}\left(
\mathfrak{Q}''^{(2)}\circ^{1\mathbf{Pth}_{\boldsymbol{\mathcal{A}}^{(2)}}}_{s} \mathfrak{Q}'^{(2)}
\right).
\tag{3}
\end{align*} 

In the just stated chain of equalities, the first equality follows from Proposition~\ref{PDPthDCatAlg}; the second equality follows by induction; finally, the last equality follows from Proposition~\ref{PDPthDCatAlg}.

Regarding the $([1],2)$-target, the following chain of equalities hold
\allowdisplaybreaks
\begin{align*}
\mathrm{tg}^{([1],2)}_{s}\left(
\mathfrak{P}''^{(2)}\circ^{1\mathbf{Pth}_{\boldsymbol{\mathcal{A}}^{(2)}}}_{s} \mathfrak{P}'^{(2)}
\right)
&=
\mathrm{tg}_{s}^{([1],2)}\left(
\mathfrak{P}''^{(2)}
\right)
\tag{1}
\\&=
\mathrm{tg}_{s}^{([1],2)}\left(
\mathfrak{Q}''^{(2)}
\right)
\tag{2}
\\&=
\mathrm{tg}^{([1],2)}_{s}\left(
\mathfrak{Q}''^{(2)}\circ^{1\mathbf{Pth}_{\boldsymbol{\mathcal{A}}^{(2)}}}_{s} \mathfrak{Q}'^{(2)}
\right).
\tag{3}
\end{align*} 

In the just stated chain of equalities, the first equality follows from Proposition~\ref{PDPthDCatAlg}; the second equality follows by induction; finally, the last equality follows from Proposition~\ref{PDPthDCatAlg}.

This completes the Case~(2.7).

This completes the Case~(2).

This completes the proof.
\end{proof}

\begin{restatable}{corollary}{CDUpsCongDZ}
\label{CDUpsCongDZ} Let $s$ be a sort in $S$ and $(\mathfrak{P}^{(2)},\mathfrak{Q}^{(2)})\in \Upsilon^{[1]}_{s}$. Then the following statements hold
\begin{enumerate}
\item[(i)] $\mathrm{sc}^{(0,2)}(\mathfrak{P}^{(2)})=\mathrm{sc}^{(0,2)}(\mathfrak{Q}^{(2)})$;
\item[(ii)] $\mathrm{tg}^{(0,2)}(\mathfrak{P}^{(2)})=\mathrm{tg}^{(0,2)}(\mathfrak{Q}^{(2)})$.
\end{enumerate}
\end{restatable}
\begin{proof}
Unravelling the mappings $\mathrm{sc}^{(0,2)}$ and $\mathrm{tg}^{(0,2)}$, introduced in Definition~\ref{DDScTgZ}, we have that 
\begin{align*}
\mathrm{sc}^{(0,2)}&=
\mathrm{sc}^{(0,[1])}
\circ
\mathrm{ip}^{([1],X)@}
\circ
\mathrm{sc}^{([1],2)};
&
\mathrm{tg}^{(0,2)}&=
\mathrm{tg}^{(0,[1])}
\circ
\mathrm{ip}^{([1],X)@}
\circ
\mathrm{tg}^{([1],2)}.
\end{align*}

Since $(\mathfrak{P}^{(2)},\mathfrak{Q}^{(2)})\in\Upsilon^{[1]}_{s}$, then, from Lemma~\ref{LDUpsCong}, we have that 
\begin{align*}
\mathrm{sc}^{([1],2)}_{s}\left(
\mathfrak{P}^{(2)}
\right)&
=\mathrm{sc}^{([1],2)}_{s}\left(
\mathfrak{Q}^{(2)}
\right);
&
\mathrm{tg}^{([1],2)}_{s}\left(
\mathfrak{P}^{(2)}
\right)
&=\mathrm{tg}^{([1],2)}_{s}\left(
\mathfrak{Q}^{(2)}
\right).
\end{align*} 

The statement follows easily.
\end{proof}

We next prove the equivalent of Corollary~\ref{CDUps} for the $\Sigma^{\boldsymbol{\mathcal{A}}^{(2)}}$-congruence $\Upsilon^{[1]}$.

\begin{restatable}{corollary}{CDUpsCongDUId}
\label{CDUpsCongDUId} Let $s$ be a sort in $S$ and $(\mathfrak{P}^{(2)},\mathfrak{Q}^{(2)})\in\Upsilon^{[1]}_{s}$. If $\mathfrak{P}^{(2)}$ or $\mathfrak{Q}^{(2)}$ is a $(2,[1])$-identity second-order path, then $\mathfrak{P}^{(2)}=\mathfrak{Q}^{(2)}$.
\end{restatable}
\begin{proof}
We recall from Remark~\ref{RDUpsCong} that 
$$\mathrm{C}^{\omega}_{\mathbf{Pth}_{\boldsymbol{\mathcal{A}}^{(2)}}}\left(
\Upsilon^{(1)}\right)
=\bigcup_{n\in\mathbb{N}}\mathrm{C}^{n}_{\mathbf{Pth}_{\boldsymbol{\mathcal{A}}^{(2)}}}\left(
\Upsilon^{(1)}
\right)
=\Upsilon^{[1]}.$$ 

We prove the statement by induction on  $n\in\mathbb{N}$.

\textsf{Base step of the induction.}

Let us recall from Remark~\ref{RDUpsCong} that 
\[
\mathrm{C}^{0}_{\mathbf{Pth}_{\boldsymbol{\mathcal{A}}^{(2)}}}\left(
\Upsilon^{(1)}\right)
=
\Upsilon^{(1)}\cup\left(\Upsilon^{(1)}\right)^{-1}\cup\Delta_{\mathrm{Pth}_{\boldsymbol{\mathcal{A}}^{(2)}}}.
\]

The statement trivially holds for pairs $(\mathfrak{P}^{(2)},\mathfrak{Q}^{(2)})\in \Delta_{\mathrm{Pth}_{\boldsymbol{\mathcal{A}}^{(2)},s}}$. 

If the pair $(\mathfrak{P}^{(2)},\mathfrak{Q}^{(2)})$ is in $\Upsilon^{(1)}_{s}$, the statement follows from Corollary~\ref{CDUps}. In case $(\mathfrak{P}^{(2)},\mathfrak{Q}^{(2)})$ is a pair in $(\Upsilon^{(1)}_{s})^{-1}$, we reason in a similar way.

This completes the base case.

\textsf{Inductive step of the induction.}

Assume that the statement holds for $n\in\mathbb{N}$, i.e., for every sort $s\in S$ and every pair of paths $(\mathfrak{P}^{(2)},\mathfrak{Q}^{(2)})$ in $\mathrm{C}^{n}_{\mathbf{Pth}_{\boldsymbol{\mathcal{A}}^{(2)}}}\left(
\Upsilon^{(1)}\right)$, if either $\mathfrak{P}^{(2)}$ or $\mathfrak{Q}^{(2)}$ is a $(2,[1])$-identity second-order path, then $\mathfrak{P}^{(2)}=\mathfrak{Q}^{(2)}$.

We now prove the statement for $n+1$. Let $s$ be a sort in $S$ and let $(\mathfrak{P}^{(2)},\mathfrak{Q}^{(2)})$ be a pair of second-order paths in $\mathrm{Pth}_{\boldsymbol{\mathcal{A}}^{(2)},s}$ such that $(\mathfrak{P}^{(2)},\mathfrak{Q}^{(2)})\in\mathrm{C}^{n+1}_{\mathbf{Pth}_{\boldsymbol{\mathcal{A}}^{(2)}}}\left(
\Upsilon^{(1)}\right)$.

Let us recall from Remark~\ref{RDUpsCong} that 
\allowdisplaybreaks
\begin{multline*}
\mathrm{C}^{n+1}_{\mathbf{Pth}_{\boldsymbol{\mathcal{A}}^{(2)}}}\left(
\Upsilon^{(1)}
\right)
=\left(
\mathrm{C}^{n}_{\mathbf{Pth}_{\boldsymbol{\mathcal{A}}^{(2)}}}\left(
\Upsilon^{(1)}
\right)
\circ
\mathrm{C}^{n}_{\mathbf{Pth}_{\boldsymbol{\mathcal{A}}^{(2)}}}\left(
\Upsilon^{(1)}
\right)
\right)
\cup
\\
\left(
\bigcup_{\gamma\in\Sigma^{\boldsymbol{\mathcal{A}}^{(2)}}_{\neq\lambda,s}}
\gamma^{\mathbf{Pth}_{\boldsymbol{\mathcal{A}}^{(2)}}}
\times
\gamma^{\mathbf{Pth}_{\boldsymbol{\mathcal{A}}^{(2)}}}
\left[
\mathrm{C}^{n}_{\mathbf{Pth}_{\boldsymbol{\mathcal{A}}^{(2)}}}\left(
\Upsilon^{(1)}
\right)_{\mathrm{ar}(\gamma)}
\right]
\right)_{s\in S}.
\end{multline*}

Then either (1),  $(\mathfrak{P}^{(2)},\mathfrak{Q}^{(2)})$ is a pair in $\mathrm{C}^{n}_{\mathbf{Pth}_{\boldsymbol{\mathcal{A}}^{(2)}}}(
\Upsilon^{(1)}
)
\circ
\mathrm{C}^{n}_{\mathbf{Pth}_{\boldsymbol{\mathcal{A}}^{(2)}}}(
\Upsilon^{(1)}
)$ or (2),  $(\mathfrak{P}^{(2)},\mathfrak{Q}^{(2)})$ is a pair in $\gamma^{\mathbf{Pth}_{\boldsymbol{\mathcal{A}}^{(2)}}}
\times
\gamma^{\mathbf{Pth}_{\boldsymbol{\mathcal{A}}^{(2)}}}
[
\mathrm{C}^{n}_{\mathbf{Pth}_{\boldsymbol{\mathcal{A}}^{(2)}}}(
\Upsilon^{(1)}
)_{\mathrm{ar}(\gamma)}]$ for some operation symbol $\gamma\in \Sigma^{\boldsymbol{\mathcal{A}}^{(2)}}_{\neq\lambda, s}$.

If (1), then there exists a second-order path $\mathfrak{R}^{(2)}\in\mathrm{Pth}_{\boldsymbol{\mathcal{A}}^{(2)},s}$ for which $(\mathfrak{P}^{(2)}, \mathfrak{R}^{(2)})$ and $(\mathfrak{R}^{(2)}, \mathfrak{Q}^{(2)})$ belong to $\mathrm{C}^{n}_{\mathbf{Pth}_{\boldsymbol{\mathcal{A}}^{(2)}}}(
\Upsilon^{(1)})$.  Since $(\mathfrak{P}^{(2)}, \mathfrak{R}^{(2)})$ is in $\mathrm{C}^{n}_{\mathbf{Pth}_{\boldsymbol{\mathcal{A}}^{(2)}}}(
\Upsilon^{(1)})_{s}$ then, by induction, we have that if either $\mathfrak{P}^{(2)}$ or $\mathfrak{R}^{(2)}$ is a $(2,[1])$-identity second-order path, then $\mathfrak{P}^{(2)}=\mathfrak{R}^{(2)}$.

Since $(\mathfrak{R}^{(2)}, \mathfrak{Q}^{(2)})$ is in $\mathrm{C}^{n}_{\mathbf{Pth}_{\boldsymbol{\mathcal{A}}^{(2)}}}(
\Upsilon^{(1)})_{s}$ then, by induction, we have that  if either $\mathfrak{R}^{(2)}$ or $\mathfrak{Q}^{(2)}$ is a $(2,[1])$-identity second-order path, then $\mathfrak{R}^{(2)}=\mathfrak{Q}^{(2)}$.

All in all, we conclude that  if either $\mathfrak{P}^{(2)}$ or $\mathfrak{Q}^{(2)}$ is a $(2,[1])$-identity second-order path, then $\mathfrak{P}^{(2)}=\mathfrak{Q}^{(2)}$.

This completes the Case~(1).

If~(2), then there exists a unique word $\mathbf{s}\in S^{\star}-\{\lambda\}$, a unique $s\in S$, a unique operation symbol $\gamma\in \Sigma^{\boldsymbol{\mathcal{A}}^{(2)}}_{\mathbf{s},s}$ and a unique family of pairs $((\mathfrak{P}^{(2)}_{j}, \mathfrak{Q}^{(2)}_{j}))_{j\in\bb{\mathbf{s}}}$ in $\mathrm{C}^{n}_{\mathbf{Pth}_{\boldsymbol{\mathcal{A}}^{(2)}}}(\Upsilon^{(1)})_{\mathbf{s}}$ for which 
\[
\left(
\mathfrak{P}^{(2)},
\mathfrak{Q}^{(2)}
\right)
=
\left(
\gamma^{\mathbf{Pth}_{\boldsymbol{\mathcal{A}}^{(2)}}}\left(\left(\mathfrak{P}^{(2)}_{j}\right)_{j\in\bb{\mathbf{s}}}\right),
\gamma^{\mathbf{Pth}_{\boldsymbol{\mathcal{A}}^{(2)}}}\left(\left(\mathfrak{Q}^{(2)}_{j}\right)_{j\in\bb{\mathbf{s}}}\right)
\right).
\]

We will distinguish the following cases according to the different possibilities for the operation symbol $\gamma\in \Sigma^{\boldsymbol{\mathcal{A}}^{(2)}}_{\mathbf{s},s}$. Note that either (2.1) $\gamma$ is an operation symbol $\sigma\in \Sigma_{\mathbf{s},s}$; or (2.2) $\gamma$ is the operation symbol of $0$-source $\mathrm{sc}^{0}_{s}\in \Sigma^{\boldsymbol{\mathcal{A}}}_{s,s}$; or (2.3) $\gamma$ is the operation symbol of $0$-target $\mathrm{tg}^{0}_{s}\in \Sigma^{\boldsymbol{\mathcal{A}}}_{s,s}$; or (2.4) $\gamma$ is the operation symbol of $0$-composition $\circ^{0}_{s}\in \Sigma^{\boldsymbol{\mathcal{A}}}_{ss,s}$; or (2.5) $\gamma$ is the operation symbol of $1$-source $\mathrm{sc}^{1}_{s}\in \Sigma^{\boldsymbol{\mathcal{A}}^{(2)}}_{s,s}$; or (2.6) $\gamma$ is the operation symbol of $1$-target $\mathrm{tg}^{1}_{s}\in \Sigma^{\boldsymbol{\mathcal{A}}^{(2)}}_{s,s}$; or (2.7) $\gamma$ is the operation symbol of $1$-composition $\circ^{1}_{s}\in \Sigma^{\boldsymbol{\mathcal{A}}^{(2)}}_{ss,s}$.

\textsf{Case~(2.1)} $\gamma$ is an operation symbol $\sigma \in \Sigma_{\mathbf{s},s}$. 

Let $\mathbf{s}$ be a word in $S^{\star}-\{\lambda\}$, let $\sigma$ be an operation symbol in $\Sigma_{\mathbf{s},s}$ and let $((\mathfrak{P}^{(2)}_{j}, \mathfrak{Q}^{(2)}_{j}))_{j\in\bb{\mathbf{s}}}$ be the family of pairs in $\mathrm{C}^{n}_{\mathbf{Pth}_{\boldsymbol{\mathcal{A}}^{(2)}}}(\Upsilon^{(1)})_{\mathbf{s}}$ for which 
\[
\left(
\mathfrak{P}^{(2)},
\mathfrak{Q}^{(2)}
\right)
=
\left(
\sigma^{\mathbf{Pth}_{\boldsymbol{\mathcal{A}}^{(2)}}}\left(\left(\mathfrak{P}^{(2)}_{j}\right)_{j\in\bb{\mathbf{s}}}\right),
\sigma^{\mathbf{Pth}_{\boldsymbol{\mathcal{A}}^{(2)}}}\left(\left(\mathfrak{Q}^{(2)}_{j}\right)_{j\in\bb{\mathbf{s}}}\right)
\right).
\]

Assume, without loss of generality, that $\mathfrak{P}^{(2)}=\sigma^{\mathbf{Pth}_{\boldsymbol{\mathcal{A}}^{(2)}}}((\mathfrak{P}^{(2)}_{j})_{j\in\bb{\mathbf{s}}})$ is a $(2,[1])$-identity second-order path. Following Claim~\ref{CDPthSigma} we have that $\bb{\mathfrak{P}^{(2)}}=\sum_{j\in\bb{\mathbf{s}}}\bb{\mathfrak{P}^{(2)}_{j}}$, then we conclude that, for every $j\in\bb{\mathbf{s}}$, it is the case that $\mathfrak{P}^{(2)}_{j}$ is a $(2,[1])$-identity second-order path. Since $((\mathfrak{P}^{(2)}_{j}, \mathfrak{Q}^{(2)}_{j}))_{j\in\bb{\mathbf{s}}}$ is a family of pairs in $\mathrm{C}^{n}_{\mathbf{Pth}_{\boldsymbol{\mathcal{A}}^{(2)}}}(\Upsilon^{(1)})_{\mathbf{s}}$, we have by induction that $\mathfrak{P}^{(2)}_{j}=\mathfrak{Q}^{(2)}_{j}$, for every $j\in\bb{\mathbf{s}}$.

The following chain of equalities holds
\allowdisplaybreaks
\begin{align*}
\mathfrak{P}^{(2)}&=
\sigma^{\mathbf{Pth}_{\boldsymbol{\mathcal{A}}^{(2)}}}\left(\left(\mathfrak{P}^{(2)}_{j}\right)_{j\in\bb{\mathbf{s}}}\right)
\tag{1}
\\&=
\sigma^{\mathbf{Pth}_{\boldsymbol{\mathcal{A}}^{(2)}}}\left(\left(\mathfrak{Q}^{(2)}_{j}\right)_{j\in\bb{\mathbf{s}}}\right)
\tag{2}
\\&=
\mathfrak{Q}^{(2)}.
\tag{3}
\end{align*}

In the just stated chain of equalities, the first equality unravels the description of $\mathfrak{P}^{(2)}$; the second equality follows from the fact that, for every $j\in\bb{\mathbf{s}}$, we have that  $\mathfrak{P}^{(2)}_{j}=\mathfrak{Q}^{(2)}_{j}$; finally, the third equality recovers the description of $\mathfrak{Q}^{(2)}$.

This completes the Case~(2.1).

\textsf{Case~(2.2)} $\gamma$ is the operation symbol $\mathrm{sc}^{0}_{s} \in \Sigma^{\boldsymbol{\mathcal{A}}}_{s,s}$. 

Let $(\mathfrak{P}'^{(2)}, \mathfrak{Q}'^{(2)})$ be the pairs in $\mathrm{C}^{n}_{\mathbf{Pth}_{\boldsymbol{\mathcal{A}}^{(2)}}}(\Upsilon^{(1)})_{s}$ for which 
\[
\left(
\mathfrak{P}^{(2)},
\mathfrak{Q}^{(2)}
\right)
=
\left(
\mathrm{sc}^{0\mathbf{Pth}_{\boldsymbol{\mathcal{A}}^{(2)}}}_{s}\left(\mathfrak{P}'^{(2)}\right),
\mathrm{sc}^{0\mathbf{Pth}_{\boldsymbol{\mathcal{A}}^{(2)}}}_{s}\left(\mathfrak{Q}'^{(2)}\right)
\right).
\]

Note that, in virtue of Claim~\ref{CDPthCatAlgScZ}, both $\mathrm{sc}^{0\mathbf{Pth}_{\boldsymbol{\mathcal{A}}^{(2)}}}_{s}(\mathfrak{P}'^{(2)})$ and $\mathrm{sc}^{0\mathbf{Pth}_{\boldsymbol{\mathcal{A}}^{(2)}}}_{s}(\mathfrak{Q}'^{(2)})$ are $(2,[1])$-identity second-order paths. The following chain of equalities hold
\allowdisplaybreaks
\begin{align*}
\mathfrak{P}^{(2)}&=
\mathrm{sc}^{0\mathbf{Pth}_{\boldsymbol{\mathcal{A}}^{(2)}}}_{s}\left(
\mathfrak{P}'^{(2)}
\right)
\tag{1}
\\&=
\mathrm{ip}^{(2,0)\sharp}_{s}\left(\mathrm{sc}^{(0,2)}_{s}\left(
\mathfrak{P}'^{(2)}
\right)\right)
\tag{2}
\\&=
\mathrm{ip}^{(2,0)\sharp}_{s}\left(\mathrm{sc}^{(0,2)}_{s}\left(
\mathfrak{Q}'^{(2)}
\right)\right)
\tag{3}
\\&=
\mathrm{sc}^{0\mathbf{Pth}_{\boldsymbol{\mathcal{A}}^{(2)}}}_{s}\left(
\mathfrak{Q}'^{(2)}
\right)
\tag{4}
\\&=
\mathfrak{Q}^{(2)}.
\end{align*}

In the just stated chain of equalities,  the first equality unravels the description of $\mathfrak{P}^{(2)}$; the second equality follows from the definition of the $0$-source operation symbol in the many-sorted partial $\Sigma^{\boldsymbol{\mathcal{A}}^{(2)}}$-algebra $\mathbf{Pth}_{\boldsymbol{\mathcal{A}}^{(2)}}$, according to Proposition~\ref{PDPthCatAlg}; the third equality follows from the fact that $(\mathfrak{P}'^{(2)},\mathfrak{Q}'^{(2)})$ are in $\Upsilon^{[1]}_{s}$. Thus, according to Corollary~\ref{CDUpsCongDZ}, they have the same $(0,2)$-source; the fourth equality recovers the definition of the $0$-source operation symbol in the many-sorted partial $\Sigma^{\boldsymbol{\mathcal{A}}^{(2)}}$-algebra $\mathbf{Pth}_{\boldsymbol{\mathcal{A}}^{(2)}}$, according to Proposition~\ref{PDPthCatAlg}; finally, the last equality recovers the description of $\mathfrak{Q}^{(2)}$.

This completes the Case~(2.2).

\textsf{Case~(2.3)} $\gamma$ is the operation symbol $\mathrm{tg}^{0}_{s} \in \Sigma^{\boldsymbol{\mathcal{A}}}_{s,s}$. 

Let $(\mathfrak{P}'^{(2)}, \mathfrak{Q}'^{(2)})$ be the pairs in $\mathrm{C}^{n}_{\mathbf{Pth}_{\boldsymbol{\mathcal{A}}^{(2)}}}(\Upsilon^{(1)})_{s}$ for which 
\[
\left(
\mathfrak{P}^{(2)},
\mathfrak{Q}^{(2)}
\right)
=
\left(
\mathrm{tg}^{0\mathbf{Pth}_{\boldsymbol{\mathcal{A}}^{(2)}}}_{s}\left(\mathfrak{P}'^{(2)}\right),
\mathrm{tg}^{0\mathbf{Pth}_{\boldsymbol{\mathcal{A}}^{(2)}}}_{s}\left(\mathfrak{Q}'^{(2)}\right)
\right).
\]

This case can be proven by a similar reasoning than that used for Case~(2.2).

This completes the Case~(2.3).

\textsf{Case~(2.4)} $\gamma$ is the operation symbol $\circ^{0}_{s} \in \Sigma^{\boldsymbol{\mathcal{A}}}_{ss,s}$. 

Let $(\mathfrak{P}'^{(2)}, \mathfrak{Q}'^{(2)})$ and $(\mathfrak{P}''^{(2)}, \mathfrak{Q}''^{(2)})$ be the pairs in $\mathrm{C}^{n}_{\mathbf{Pth}_{\boldsymbol{\mathcal{A}}^{(2)}}}(\Upsilon^{(1)})_{s}$ satisfying 
\begin{align*}
\mathrm{sc}^{(0,2)}_{s}\left(
\mathfrak{P}''^{(2)}
\right)
&=
\mathrm{tg}^{(0,2)}_{s}\left(
\mathfrak{P}'^{(2)}
\right);
&
\mathrm{sc}^{(0,2)}_{s}\left(
\mathfrak{Q}''^{(2)}
\right)
&=
\mathrm{tg}^{(0,2)}_{s}\left(
\mathfrak{Q}'^{(2)}
\right),
\end{align*}
for which 
\[
\left(
\mathfrak{P}^{(2)},
\mathfrak{Q}^{(2)}
\right)
=
\left(
\mathfrak{P}''^{(2)}\circ^{0\mathbf{Pth}_{\boldsymbol{\mathcal{A}}^{(2)}}}_{s} \mathfrak{P}'^{(2)},
\mathfrak{Q}''^{(2)}\circ^{0\mathbf{Pth}_{\boldsymbol{\mathcal{A}}^{(2)}}}_{s} \mathfrak{Q}'^{(2)}
\right).
\]

Assume, without loss of generality, that $\mathfrak{P}^{(2)}=\mathfrak{P}''^{(2)}\circ^{0\mathbf{Pth}_{\boldsymbol{\mathcal{A}}^{(2)}}}_{s} \mathfrak{P}'^{(2)}$ is a $(2,[1])$-identity second-order path. Following Claim~\ref{CDPthCatAlgCompZ} we have that $\bb{\mathfrak{P}^{(2)}}=\bb{\mathfrak{P}'^{(2)}}+\bb{\mathfrak{P}''^{(2)}}$, then we conclude that both $\mathfrak{P}'^{(2)}$ and $\mathfrak{P}''^{(2)}$ are $(2,[1])$-identity second-order paths. Since both $(\mathfrak{P}'^{(2)}, \mathfrak{Q}'^{(2)})$ and $(\mathfrak{P}''^{(2)}, \mathfrak{Q}''^{(2)})$ are pairs in $\mathrm{C}^{n}_{\mathbf{Pth}_{\boldsymbol{\mathcal{A}}^{(2)}}}(\Upsilon^{(1)})_{s}$, we have by induction that $\mathfrak{P}'^{(2)}=\mathfrak{Q}'^{(2)}$ and $\mathfrak{P}''^{(2)}=\mathfrak{Q}''^{(2)}$.

The following chain of equalities holds
\allowdisplaybreaks
\begin{align*}
\mathfrak{P}^{(2)}&=
\mathfrak{P}''^{(2)}\circ^{0\mathbf{Pth}_{\boldsymbol{\mathcal{A}}^{(2)}}}_{s} \mathfrak{P}'^{(2)}
\tag{1}
\\&=
\mathfrak{Q}''^{(2)}\circ^{0\mathbf{Pth}_{\boldsymbol{\mathcal{A}}^{(2)}}}_{s} \mathfrak{Q}'^{(2)}
\tag{2}
\\&=
\mathfrak{Q}^{(2)}.
\tag{3}
\end{align*}

In the just stated chain of equalities, the first equality unravels the description of $\mathfrak{P}^{(2)}$; the second equality follows from the fact that $\mathfrak{P}'^{(2)}=\mathfrak{Q}'^{(2)}$ and $\mathfrak{P}''^{(2)}=\mathfrak{Q}''^{(2)}$; finally, the third equality recovers the description of $\mathfrak{Q}^{(2)}$.

This completes the Case~(2.4).

\textsf{Case~(2.5)} $\gamma$ is the operation symbol $\mathrm{sc}^{1}_{s} \in \Sigma^{\boldsymbol{\mathcal{A}}}_{s,s}$. 

Let $(\mathfrak{P}'^{(2)}, \mathfrak{Q}'^{(2)})$ be the pairs in $\mathrm{C}^{n}_{\mathbf{Pth}_{\boldsymbol{\mathcal{A}}^{(2)}}}(\Upsilon^{(1)})_{s}$ for which 
\[
\left(
\mathfrak{P}^{(2)},
\mathfrak{Q}^{(2)}
\right)
=
\left(
\mathrm{sc}^{1\mathbf{Pth}_{\boldsymbol{\mathcal{A}}^{(2)}}}_{s}\left(\mathfrak{P}'^{(2)}\right),
\mathrm{sc}^{1\mathbf{Pth}_{\boldsymbol{\mathcal{A}}^{(2)}}}_{s}\left(\mathfrak{Q}'^{(2)}\right)
\right).
\]

Note that, in virtue of Proposition~\ref{PDPthDCatAlg}, both $\mathrm{sc}^{1\mathbf{Pth}_{\boldsymbol{\mathcal{A}}^{(2)}}}_{s}(\mathfrak{P}'^{(2)})$ and $\mathrm{sc}^{1\mathbf{Pth}_{\boldsymbol{\mathcal{A}}^{(2)}}}_{s}(\mathfrak{Q}'^{(2)})$ are $(2,[1])$-identity second-order paths. The following chain of equalities hold
\allowdisplaybreaks
\begin{align*}
\mathfrak{P}^{(2)}&=
\mathrm{sc}^{1\mathbf{Pth}_{\boldsymbol{\mathcal{A}}^{(2)}}}_{s}\left(
\mathfrak{P}'^{(2)}
\right)
\tag{1}
\\&=
\mathrm{ip}^{(2,[1])\sharp}_{s}\left(\mathrm{sc}^{([1],2)}_{s}\left(
\mathfrak{P}'^{(2)}
\right)\right)
\tag{2}
\\&=
\mathrm{ip}^{(2,[1])\sharp}_{s}\left(\mathrm{sc}^{([1],2)}_{s}\left(
\mathfrak{Q}'^{(2)}
\right)\right)
\tag{3}
\\&=
\mathrm{sc}^{1\mathbf{Pth}_{\boldsymbol{\mathcal{A}}^{(2)}}}_{s}\left(
\mathfrak{Q}'^{(2)}
\right)
\tag{4}
\\&=
\mathfrak{Q}^{(2)}.
\end{align*}

In the just stated chain of equalities,  the first equality unravels the description of $\mathfrak{P}^{(2)}$; the second equality follows from the definition of the $1$-source operation symbol in the many-sorted partial $\Sigma^{\boldsymbol{\mathcal{A}}^{(2)}}$-algebra $\mathbf{Pth}_{\boldsymbol{\mathcal{A}}^{(2)}}$, according to Proposition~\ref{PDPthDCatAlg}; the third equality follows from the fact that $(\mathfrak{P}'^{(2)},\mathfrak{Q}'^{(2)})$ are in $\Upsilon^{[1]}_{s}$. Thus, according to Lemma~\ref{LDUpsCong}, they have the same $([1],2)$-source; the fourth equality recovers the definition of the $1$-source operation symbol in the many-sorted partial $\Sigma^{\boldsymbol{\mathcal{A}}^{(2)}}$-algebra $\mathbf{Pth}_{\boldsymbol{\mathcal{A}}^{(2)}}$, according to Proposition~\ref{PDPthDCatAlg}; finally, the last equality recovers the description of $\mathfrak{Q}^{(2)}$.

This completes the Case~(2.5).

\textsf{Case~(2.6)} $\gamma$ is the operation symbol $\mathrm{tg}^{1}_{s} \in \Sigma^{\boldsymbol{\mathcal{A}}}_{s,s}$. 

Let $(\mathfrak{P}'^{(2)}, \mathfrak{Q}'^{(2)})$ be the pairs in $\mathrm{C}^{n}_{\mathbf{Pth}_{\boldsymbol{\mathcal{A}}^{(2)}}}(\Upsilon^{(1)})_{s}$ for which 
\[
\left(
\mathfrak{P}^{(2)},
\mathfrak{Q}^{(2)}
\right)
=
\left(
\mathrm{tg}^{1\mathbf{Pth}_{\boldsymbol{\mathcal{A}}^{(2)}}}_{s}\left(\mathfrak{P}'^{(2)}\right),
\mathrm{tg}^{1\mathbf{Pth}_{\boldsymbol{\mathcal{A}}^{(2)}}}_{s}\left(\mathfrak{Q}'^{(2)}\right)
\right).
\]

This case can be proven by a similar reasoning than that used for Case~(2.5).

This completes the Case~(2.6).

\textsf{Case~(2.7)} $\gamma$ is the operation symbol $\circ^{1}_{s} \in \Sigma^{\boldsymbol{\mathcal{A}}}_{ss,s}$. 

Let $(\mathfrak{P}'^{(2)}, \mathfrak{Q}'^{(2)})$ and $(\mathfrak{P}''^{(2)}, \mathfrak{Q}''^{(2)})$ be the pairs in $\mathrm{C}^{n}_{\mathbf{Pth}_{\boldsymbol{\mathcal{A}}^{(2)}}}(\Upsilon^{(1)})_{s}$ satisfying 
\begin{align*}
\mathrm{sc}^{(0,2)}_{s}\left(
\mathfrak{P}''^{(2)}
\right)
&=
\mathrm{tg}^{(0,2)}_{s}\left(
\mathfrak{P}'^{(2)}
\right);
&
\mathrm{sc}^{(0,2)}_{s}\left(
\mathfrak{Q}''^{(2)}
\right)
&=
\mathrm{tg}^{(0,2)}_{s}\left(
\mathfrak{Q}'^{(2)}
\right),
\end{align*}
for which 
\[
\left(
\mathfrak{P}^{(2)},
\mathfrak{Q}^{(2)}
\right)
=
\left(
\mathfrak{P}''^{(2)}\circ^{1\mathbf{Pth}_{\boldsymbol{\mathcal{A}}^{(2)}}}_{s} \mathfrak{P}'^{(2)},
\mathfrak{Q}''^{(2)}\circ^{1\mathbf{Pth}_{\boldsymbol{\mathcal{A}}^{(2)}}}_{s} \mathfrak{Q}'^{(2)}
\right).
\]

Assume, without loss of generality, that $\mathfrak{P}^{(2)}=\mathfrak{P}''^{(2)}\circ^{1\mathbf{Pth}_{\boldsymbol{\mathcal{A}}^{(2)}}}_{s} \mathfrak{P}'^{(2)}$ is a $(2,[1])$-identity second-order path. Following Proposition~\ref{PDPthDCatAlg}, we have that $\bb{\mathfrak{P}^{(2)}}=\bb{\mathfrak{P}'^{(2)}}+\bb{\mathfrak{P}''^{(2)}}$, then we conclude that both $\mathfrak{P}'^{(2)}$ and $\mathfrak{P}''^{(2)}$ are $(2,[1])$-identity second-order paths. Since both $(\mathfrak{P}'^{(2)}, \mathfrak{Q}'^{(2)})$ and $(\mathfrak{P}''^{(2)}, \mathfrak{Q}''^{(2)})$ are pairs in $\mathrm{C}^{n}_{\mathbf{Pth}_{\boldsymbol{\mathcal{A}}^{(2)}}}(\Upsilon^{(1)})_{s}$, we have by induction that $\mathfrak{P}'^{(2)}=\mathfrak{Q}'^{(2)}$ and $\mathfrak{P}''^{(2)}=\mathfrak{Q}''^{(2)}$.

The following chain of equalities holds
\allowdisplaybreaks
\begin{align*}
\mathfrak{P}^{(2)}&=
\mathfrak{P}''^{(2)}\circ^{1\mathbf{Pth}_{\boldsymbol{\mathcal{A}}^{(2)}}}_{s} \mathfrak{P}'^{(2)}
\tag{1}
\\&=
\mathfrak{Q}''^{(2)}\circ^{1\mathbf{Pth}_{\boldsymbol{\mathcal{A}}^{(2)}}}_{s} \mathfrak{Q}'^{(2)}
\tag{2}
\\&=
\mathfrak{Q}^{(2)}.
\tag{3}
\end{align*}

In the just stated chain of equalities, the first equality unravels the description of $\mathfrak{P}^{(2)}$; the second equality follows from the fact that $\mathfrak{P}'^{(2)}=\mathfrak{Q}'^{(2)}$ and $\mathfrak{P}''^{(2)}=\mathfrak{Q}''^{(2)}$; finally, the third equality recovers the description of $\mathfrak{Q}^{(2)}$.

This completes the Case~(2.7).

This completes the Case~(2).

This completes the proof.
\end{proof}

\section{
\texorpdfstring
{A structure of partial $\Sigma^{\boldsymbol{\mathcal{A}}^{(2)}}$-algebra on $[\mathrm{Pth}_{\boldsymbol{\mathcal{A}}^{(2)}}]_{\Upsilon^{[1]}}$}
{A partial algebra on the Upsilon quotient}
}
In this section we only acquaint a consequence of the main result of the previous section. Since 
$\Upsilon^{[1]}$ is a $\Sigma^{\boldsymbol{\mathcal{A}}^{(2)}}$-congruence on $\mathbf{Pth}_{\boldsymbol{\mathcal{A}}^{(2)}}$ we can endow the many-sorted quotient $\mathrm{Pth}_{\boldsymbol{\mathcal{A}}^{(2)}}/{\Upsilon^{[1]}}$, that we will simply denote by $[\mathrm{Pth}_{\boldsymbol{\mathcal{A}}^{(2)}}]_{\Upsilon^{[1]}}$, with a structure of many-sorted partial $\Sigma^{\boldsymbol{\mathcal{A}}^{(2)}}$-algebra.
\index{path!second-order!$[\mathrm{Pth}_{\boldsymbol{\mathcal{A}}^{(2)}}]_{\Upsilon^{[1]}}$}

\begin{restatable}{proposition}{PDUpsCongCatAlg}
\label{PDUpsCongCatAlg}
\index{path!second-order!$[\mathbf{Pth}_{\boldsymbol{\mathcal{A}}^{(2)}}]_{\Upsilon^{[1]}}$}
 The $S$-sorted set $[\mathrm{Pth}_{\boldsymbol{\mathcal{A}}^{(2)}}]_{\Upsilon^{[1]}}$ is equipped, in a natural way, with a structure of many-sorted partial $\Sigma^{\boldsymbol{\mathcal{A}}^{(2)}}$-algebra. We denote by $[\mathbf{Pth}_{\boldsymbol{\mathcal{A}}^{(2)}}]_{\Upsilon^{[1]}}$ the corresponding $\Sigma^{\boldsymbol{\mathcal{A}}^{(2)}}$-algebra (which is a quotient of $\mathbf{Pth}_{\boldsymbol{\mathcal{A}}^{(2)}}$, the many-sorted partial $\Sigma^{\boldsymbol{\mathcal{A}}^{(2)}}$-algebra constructed in Proposition~\ref{PDPthDCatAlg}).

\index{projection!second-order!$\mathrm{pr}^{\Upsilon^{[1]}}$}
The mapping $\mathrm{pr}^{\Upsilon^{[1]}}$, that is,
\[
\mathrm{pr}^{\Upsilon^{[1]}}
\colon
\mathbf{Pth}_{\boldsymbol{\mathcal{A}}^{(2)}}
\longrightarrow
[\mathbf{Pth}_{\boldsymbol{\mathcal{A}}^{(2)}}]_{\Upsilon^{[1]}}
\] 
is a closed and surjective $\Sigma^{\boldsymbol{\mathcal{A}}^{(2)}}$-homomorphism from $\mathbf{Pth}_{\boldsymbol{\mathcal{A}}^{(2)}}$ to $[\mathbf{Pth}_{\boldsymbol{\mathcal{A}}^{(2)}}]_{\Upsilon^{[1]}}$.
\end{restatable}

\section{
\texorpdfstring
{An Artinian preorder on $\coprod[\mathrm{Pth}_{\boldsymbol{\mathcal{A}}^{(2)}}]_{\Upsilon^{[1]}}$}
{An Artinian preorder on the Upsilon congruence}
}

The aim of this section is to provide the coproduct of the many-sorted set of second-order path classes $\coprod[\mathrm{Pth}_{\boldsymbol{\mathcal{A}}^{(2)}}]_{\Upsilon^{[1]}}$ with an Artinian preorder. 

\begin{restatable}{definition}{DDUpsOrd}
\label{DDUpsOrd}
\index{partial preorder!second-order!$\leq_{[\mathbf{Pth}_{\boldsymbol{\mathcal{A}}^{(2)}}]_{\Upsilon^{[1]}}}$}
Let $\leq_{[\mathbf{Pth}_{\boldsymbol{\mathcal{A}}^{(2)}}]_{\Upsilon^{[1]}}}$ be the binary relation defined on $\coprod[\mathrm{Pth}_{\boldsymbol{\mathcal{A}}^{(2)}}]_{\Upsilon^{[1]}}$ containing every pair $(([\mathfrak{Q}^{(2)}]_{\Upsilon^{[1]}_{t}},t), ([\mathfrak{P}^{(2)}]_{\Upsilon^{[1]}_{s}},s))$ in $(\coprod[\mathrm{Pth}_{\boldsymbol{\mathcal{A}}^{(2)}}]_{\Upsilon^{[1]}})^{2}$ satisfying that there exists a natural number $m\in\mathbb{N}-\{0\}$, a word $\mathbf{w}\in S^{\star}$ of length $\bb{\mathbf{w}}=m+1$ and a family of second-order paths $(\mathfrak{R}^{(2)}_{k})_{k\in\bb{\mathbf{w}}}$ in $\mathrm{Pth}_{\boldsymbol{\mathcal{A}}^{(2)},\mathbf{w}}$ such that $w_{0}=t$, $[\mathfrak{R}^{(2)}_{0}]_{\Upsilon^{[1]}_{t}}=[\mathfrak{Q}^{(2)}]_{\Upsilon^{[1]}_{t}}$, $w_{m}=s$, $[\mathfrak{R}^{(2)}_{m}]_{\Upsilon^{[1]}_{s}}=[\mathfrak{P}^{(2)}]_{\Upsilon^{[1]}_{s}}$ and, for every $k\in m$, $w_{k}=w_{k+1}$ and $(\mathfrak{R}^{(2)}_{k}, \mathfrak{R}^{(2)}_{k+1})\in \Upsilon^{[1]}_{w_{k}}$ or $(\mathfrak{R}^{(2)}_{k}, w_{k})\leq_{\mathbf{Pth}_{\boldsymbol{\mathcal{A}}^{(2)}}} (\mathfrak{R}^{(2)}_{k+1}, w_{k+1})$.
\end{restatable}

The aim of this section is to prove that $\leq_{[\mathbf{Pth}_{\boldsymbol{\mathcal{A}}^{(2)}}]_{\Upsilon^{[1]}}}$ is an Artinian preorder. This is what we do in the following Proposition.

\begin{restatable}{proposition}{PDUpsOrdArt}
\label{PDUpsOrdArt} $(\coprod[\mathrm{Pth}_{\boldsymbol{\mathcal{A}}^{(2)}}]_{\Upsilon^{[1]}}, \leq_{[\mathbf{Pth}_{\boldsymbol{\mathcal{A}}^{(2)}}]_{\Upsilon^{[1]}}})$ is a partially preordered set. Moreover, in this partially preordered set there is not any strictly decreasing $\omega_{0}$-chain, i.e., $(\coprod[\mathrm{Pth}_{\boldsymbol{\mathcal{A}}^{(2)}}]_{\Upsilon^{[1]}}, \leq_{[\mathbf{Pth}_{\boldsymbol{\mathcal{A}}^{(2)}}]_{\Upsilon^{[1]}}})$ is an Artinian preordered set.
\end{restatable}
\begin{proof}
That $\leq_{[\mathbf{Pth}_{\boldsymbol{\mathcal{A}}^{(2)}}]_{\Upsilon^{[1]}}}$ is reflexive follows from the fact that either $\Upsilon^{[1]}$ or $\leq_{\mathbf{Pth}_{\boldsymbol{\mathcal{A}}^{(2)}}}$ are reflexive. Indeed, let $s$ be a sort in $S$ and let $[\mathfrak{P}^{(2)}]_{\Upsilon^{[1]}_{s}}$ be a second-order path class in $[\mathrm{Pth}_{\boldsymbol{\mathcal{A}}^{(2)},s}]_{\Upsilon^{[1]}_{s}}$. Then, for the sequence $(\mathfrak{P}^{(2)},\mathfrak{P}^{(2)})$ in $\mathrm{Pth}_{\boldsymbol{\mathcal{A}}^{(2)},ss}$, we have that 
\[
\left(
\left[\mathfrak{P}^{(2)}\right]_{\Upsilon^{[1]}_{s}}, s\right)
\leq_{[\mathbf{Pth}_{\boldsymbol{\mathcal{A}}^{(2)}}]_{\Upsilon^{[1]}}} 
\left(
\left[\mathfrak{P}^{(2)}\right]_{\Upsilon^{[1]}_{s}}, s\right).
\]

That $\leq_{[\mathbf{Pth}_{\boldsymbol{\mathcal{A}}^{(2)}}]_{\Upsilon^{[1]}}}$ is transitive follows as well. Indeed, let $s,t,u$ be sorts in $S$ and let $[\mathfrak{P}^{(2)}]_{\Upsilon^{[1]}_{s}}$, $[\mathfrak{Q}^{(2)}]_{\Upsilon^{[1]}_{t}}$ and $[\mathfrak{S}^{(2)}]_{\Upsilon^{[1]}_{u}}$ be second-order path classes satisfying that 
\[
\left(\left[
\mathfrak{Q}^{(2)}
\right]_{\Upsilon^{[1]}_{t}},t\right)
\leq_{[\mathbf{Pth}_{\boldsymbol{\mathcal{A}}^{(2)}}]_{\Upsilon^{[1]}}} 
\left(\left[
\mathfrak{P}^{(2)}
\right]_{\Upsilon^{[1]}_{s}},s\right)
\]
\[
\left(\left[
\mathfrak{S}^{(2)}
\right]_{\Upsilon^{[1]}_{u}},u\right)
\leq_{[\mathbf{Pth}_{\boldsymbol{\mathcal{A}}^{(2)}}]_{\Upsilon^{[1]}}} 
\left(\left[
\mathfrak{Q}^{(2)}
\right]_{\Upsilon^{[1]}_{t}},t\right)
\]
Thus there exists a natural number $m\in\mathbb{N}-\{0\}$, a word $\mathbf{w}\in S^{\star}$ of length $\bb{\mathbf{w}}=m+1$ and a family of second-order paths $(\mathfrak{R}^{(2)}_{k})_{k\in\bb{\mathbf{w}}}$ in $\mathrm{Pth}_{\boldsymbol{\mathcal{A}}^{(2)},\mathbf{w}}$ such that $w_{0}=t$, $[\mathfrak{R}^{(2)}_{0}]_{\Upsilon^{[1]}_{t}}=[\mathfrak{Q}^{(2)}]_{\Upsilon^{[1]}_{t}}$, $w_{m}=s$, $[\mathfrak{R}^{(2)}_{m}]_{\Upsilon^{[1]}_{s}}=[\mathfrak{P}^{(2)}]_{\Upsilon^{[1]}_{s}}$ and, for every $k\in m$, $w_{k}=w_{k+1}$ and $(\mathfrak{R}^{(2)}_{k}, \mathfrak{R}^{(2)}_{k+1})\in \Upsilon^{[1]}_{w_{k}}$ or $(\mathfrak{R}^{(2)}_{k}, w_{k})\leq_{\mathbf{Pth}_{\boldsymbol{\mathcal{A}}^{(2)}}} (\mathfrak{R}^{(2)}_{k+1}, w_{k+1})$. Regarding the other inequality, there exists a natural number $m'\in\mathbb{N}-\{0\}$, a word $\mathbf{w}'\in S^{\star}$ of length $\bb{\mathbf{w}'}=m'+1$ and a family of second-order paths $(\mathfrak{R}'^{(2)}_{k})_{k\in\bb{\mathbf{w}'}}$ in $\mathrm{Pth}_{\boldsymbol{\mathcal{A}}^{(2)},\mathbf{w}'}$ such that $w'_{0}=u$, $[\mathfrak{R}'^{(2)}_{0}]_{\Upsilon^{[1]}_{u}}=[\mathfrak{S}^{(2)}]_{\Upsilon^{[1]}_{u}}$, $w'_{m'}=t$, $[\mathfrak{R}'^{(2)}_{m'}]_{\Upsilon^{[1]}_{t}}=[\mathfrak{Q}^{(2)}]_{\Upsilon^{[1]}_{t}}$ and, for every $k\in m'$, $w'_{k}=w'_{k+1}$ and $(\mathfrak{R}'^{(2)}_{k}, \mathfrak{R}'^{(2)}_{k+1})\in \Upsilon^{[1]}_{w'_{k}}$ or $(\mathfrak{R}'^{(2)}_{k}, w'_{k})\leq_{\mathbf{Pth}_{\boldsymbol{\mathcal{A}}^{(2)}}} (\mathfrak{R}'^{(2)}_{k+1}, w'_{k+1})$. 

Thus, the concatenation of the sequences $(\mathfrak{R}^{(2)}_{k})_{k\in\bb{\mathbf{w}}}$ and $(\mathfrak{R}'^{(2)}_{k})_{k\in\bb{\mathbf{w}'}}$ instantiates that 
\[
\left(\left[
\mathfrak{S}^{(2)}
\right]_{\Upsilon^{[1]}_{u}},u\right)
\leq_{[\mathbf{Pth}_{\boldsymbol{\mathcal{A}}^{(2)}}]_{\Upsilon^{[1]}}} 
\left(\left[
\mathfrak{P}^{(2)}
\right]_{\Upsilon^{[1]}_{s}},s\right).
\]

This proves that 
$\leq_{[\mathbf{Pth}_{\boldsymbol{\mathcal{A}}^{(2)}}]_{\Upsilon^{[1]}}}$ is a preorder on  $\coprod[\mathrm{Pth}_{\boldsymbol{\mathcal{A}}^{(2)}}]_{\Upsilon^{[1]}}$.

That it is Artinian follows from the fact that any strictly decreasing chain in $(\coprod[\mathrm{Pth}_{\boldsymbol{\mathcal{A}}^{(2)}}]_{\Upsilon^{[1]}}, \leq_{[\mathbf{Pth}_{\boldsymbol{\mathcal{A}}^{(2)}}]_{\Upsilon^{[1]}}})$ is given by a strictly decreasing chain that only uses the strict order $<_{\mathbf{Pth}_{\boldsymbol{\mathcal{A}}^{(2)}}}$ at each step. The statement follows from the fact that, by Proposition~\ref{PDOrdArt}, $(\coprod\mathrm{Pth}_{\boldsymbol{\mathcal{A}}^{(2)}}, \leq_{\mathbf{Pth}_{\boldsymbol{\mathcal{A}}^{(2)}}})$ is an Artinian preordered set.

This completes the proof.
\end{proof}

Taking into account the definitions above, we now prove that the coproduct of the projection mapping becomes order-preserving.

\begin{restatable}{proposition}{PDUpsProjMono}
\label{PDUpsProjMono} The mapping $\coprod\mathrm{pr}^{\Upsilon^{[1]}}$ is monotone
\[
\textstyle
\coprod\mathrm{pr}^{\Upsilon^{[1]}}\colon 
\left(
\coprod \mathrm{Pth}_{\boldsymbol{\mathcal{A}}^{(2)}},
\leq_{\mathbf{Pth}_{\boldsymbol{\mathcal{A}}^{(2)}}}
\right)
\mor
\left(
\coprod [\mathrm{Pth}_{\boldsymbol{\mathcal{A}}^{(2)}}]_{\Upsilon^{[1]}},
\leq_{[\mathbf{Pth}_{\boldsymbol{\mathcal{A}}^{(2)}}]_{\Upsilon^{[1]}}}
\right)
\]
\end{restatable}
\begin{proof}
Let $s$ and $t$ be sorts in $S$ and let $(\mathfrak{Q}^{(2)},t)$ and $(\mathfrak{P}^{(2)},s)$ be pairs in $\coprod \mathrm{Pth}_{\boldsymbol{\mathcal{A}}^{(2)}}$ satisfying that 
\[
\left(
\mathfrak{Q}^{(2)},t
\right)
\leq_{\mathbf{Pth}_{\boldsymbol{\mathcal{A}}^{(2)}}}
\left(
\mathfrak{P}^{(2)},s
\right),
\]
then the sequence $(\mathfrak{Q}^{(2)},\mathfrak{P}^{(2)})$ can be used to instantiate that 
\[
\left(
[\mathfrak{Q}^{(2)}]_{\Upsilon^{[1]}_{t}},t
\right)
\leq_{[\mathbf{Pth}_{\boldsymbol{\mathcal{A}}^{(2)}}]_{\Upsilon^{[1]}}}
\left(
[\mathfrak{P}^{(2)}]_{\Upsilon^{[1]}_{s}},s
\right).
\]

It follows that the mapping $\coprod\mathrm{pr}^{\Upsilon^{[1]}}$ is monotone.
\end{proof}

\chapter{
\texorpdfstring
{On the quotient $\llbracket \mathrm{Pth}_{\boldsymbol{\mathcal{A}}^{(2)}}\rrbracket$}
{On the second-order quotient}
}\label{S2H}

In this chapter we study $\mathrm{Ker}(\mathrm{CH}^{(2)})\vee \Upsilon^{[1]}$, the supremum of the two $\Sigma^{\boldsymbol{\mathcal{A}}^{(2)}}$-congruences, $\mathrm{Ker}(\mathrm{CH}^{(2)})$ and $\Upsilon^{[1]}$, defined on $\mathbf{Pth}_{\boldsymbol{\mathcal{A}}^{(2)}}$ the many-sorted partial $\Sigma^{\boldsymbol{\mathcal{A}}^{(2)}}$-algebra of second-order paths. For simplicity, the quotient will be denoted by $\llbracket \mathrm{Pth}_{\boldsymbol{\mathcal{A}}^{(2)}}\rrbracket$. Following this simplification the equivalence class of a second-order path $\mathfrak{P}^{(2)}\in\mathrm{Pth}_{\boldsymbol{\mathcal{A}}^{(2)},s}$, will simply be denoted by $\llbracket \mathfrak{P}^{(2)}\rrbracket_{s}$. We first prove that pairs of second-order paths in $\mathrm{Ker}(\mathrm{CH}^{(2)})\vee \Upsilon^{[1]}$ have the same length, the same $([1],2)$-source and $([1],2)$-target. As a result, pairs of second-order paths in $\mathrm{Ker}(\mathrm{CH}^{(2)})\vee \Upsilon^{[1]}$ also have the same $(0,2)$-source and $(0,2)$-target. We will also check that $(2,[1])$-identity second-order paths in $\mathrm{Ker}(\mathrm{CH}^{(2)})\vee \Upsilon^{[1]}$ are equal. Since ${\mathrm{Ker}(\mathrm{CH}^{(2)})}$ and $\Upsilon^{[1]}$ are closed $\Sigma^{\boldsymbol{\mathcal{A}}^{(2)}}$-congruence then so is $\mathrm{Ker}(\mathrm{CH}^{(2)})\vee \Upsilon^{[1]}$. Thus, we can consider the natural structure of many-sorted partial $\Sigma^{\boldsymbol{\mathcal{A}}^{(2)}}$-congruence on the quotient many-sorted set $\llbracket \mathrm{Pth}_{\boldsymbol{\mathcal{A}}^{(2)}}\rrbracket$, that we will denote by $\llbracket\mathbf{Pth}_{\boldsymbol{\mathcal{A}}^{(2)}}\rrbracket$. We then introduce the operations of $([1],2)$-source, $([1],2)$-target and $([1],2)$-identity second-order path in the quotient, that will be denoted, respectively, by $\mathrm{sc}^{([1],\llbracket 2\rrbracket)}$, $\mathrm{tg}^{([1],\llbracket 2\rrbracket)} $ and $\mathrm{ip}^{(\llbracket 2\rrbracket,[1])\sharp}$. As it was shown before, these mappings are also $\Sigma^{\boldsymbol{\mathcal{A}}}$-homomorphisms. We then introduce the operations of $(0,2)$-source, $(0,2)$-target and $(0,2)$-identity second-order path in the quotient, that will be denoted, respectively, by $\mathrm{sc}^{(0,\llbracket 2\rrbracket)}$, $\mathrm{tg}^{(0,\llbracket 2\rrbracket)} $ and $\mathrm{ip}^{(\llbracket 2\rrbracket,0)\sharp}$. As it was shown before, these mappings are also $\Sigma$-homomorphisms. The quotient $\llbracket \mathrm{Pth}_{\boldsymbol{\mathcal{A}}^{(2)}} \rrbracket$ has more structure than the original many-sorted set of second-order paths. In fact, we prove that $\llbracket \mathrm{Pth}_{\boldsymbol{\mathcal{A}}^{(2)}} \rrbracket$ has an structure of $S$-sorted $2$-category. Indeed, we prove that it is a $2$-categorial $\Sigma$-algebra, that we will denote by $\llbracket \mathsf{Pth}_{\boldsymbol{\mathcal{A}}^{(2)}} \rrbracket$. We conclude this chapter by investigating how to order the quotient $\llbracket \mathrm{Pth}_{\boldsymbol{\mathcal{A}}^{(2)}} \rrbracket$. We first introduce the relation $\leq_{\llbracket \mathbf{Pth}_{\boldsymbol{\mathcal{A}}^{(2)}} \rrbracket}$. We verify that $\leq_{\llbracket \mathbf{Pth}_{\boldsymbol{\mathcal{A}}^{(2)}} \rrbracket}$ is itself an Artinian partial preorder in the quotient. In addition, the mapping $\mathrm{pr}^{\llbracket \cdot \rrbracket}$, of projection to the class, is order-preserving.


\begin{restatable}{convention}{CDVPthConv}
\label{CDVPthConv} 
\index{path!second-order!$\llbracket\mathfrak{P}^{(2)}\rrbracket_{s}$}
\index{path!second-order!$\llbracket \mathrm{Pth}_{\boldsymbol{\mathcal{A}}^{(2)}}\rrbracket$}
In order to simplify the presentation, for a sort $s\in S$ and a second-order path $\mathfrak{P}^{(2)}\in\mathrm{Pth}_{\boldsymbol{\mathcal{A}}^{(2)},s}$, we will let $\llbracket      \mathfrak{P}^{(2)}\rrbracket_{s}$ stand for $[\mathfrak{P}^{(2)}]_{\mathrm{Ker}(\mathrm{CH}^{(2)})_{s}\vee \Upsilon^{[1]}_{s}}$, the $\mathrm{Ker}(\mathrm{CH}^{(2)})_{s}\vee \Upsilon^{[1]}_{s}$-class of $\mathfrak{P}^{(2)}$, and we will denote by $\llbracket \mathrm{Pth}_{\boldsymbol{\mathcal{A}}^{(2)}}\rrbracket$ the many-sorted quotient
\[
\llbracket \mathrm{Pth}_{\boldsymbol{\mathcal{A}}^{(2)}}\rrbracket=
\mathrm{Pth}_{\boldsymbol{\mathcal{A}}^{(2)}}/{\mathrm{Ker}(\mathrm{CH}^{(2)})\vee \Upsilon^{[1]}}.
\]
\end{restatable}

In the following lemma we prove that two second-order paths of the same sort that are related under the supremum of the Kernel of the second-order Curry-Howard mapping and the congruence $\Upsilon^{[1]}$ must have the same length, the same $([1],2)$-source and the same $([1],2)$-target.

\begin{restatable}{lemma}{LDV}
\label{LDV} Let $s$ be a sort in $S$ and $(\mathfrak{P}^{(2)},\mathfrak{Q}^{(2)})\in \mathrm{Ker}(\mathrm{CH}^{(2)}_{s})\vee \Upsilon^{[1]}_{s}$. The following statements hold.
\begin{enumerate}
\item[(i)] $\bb{\mathfrak{P}^{(2)}}=\bb{\mathfrak{Q}^{(2)}}$;
\item[(ii)] $\mathrm{sc}_{s}^{([1],2)}(\mathfrak{P}^{(2)})=\mathrm{sc}_{s}^{([1],2)}(\mathfrak{Q}^{(2)})$;
\item[(iii)] $\mathrm{tg}_{s}^{([1],2)}(\mathfrak{P}^{(2)})=\mathrm{tg}_{s}^{([1],2)}(\mathfrak{Q}^{(2)})$.
\end{enumerate}
\end{restatable}
\begin{proof}
Let us recall that, for a sort $s\in S$, two second-order paths $\mathfrak{P}^{(2)}$ and $\mathfrak{Q}^{(2)}$ in $\mathrm{Pth}_{\boldsymbol{\mathcal{A}}^{(2)},s}$ are related under $\mathrm{Ker}(\mathrm{CH}^{(2)}_{s})\vee \Upsilon^{[1]}_{s}$ if, and only if, there exists a natural number $m\in \mathbb{N}-\{0\}$ and a sequence of second-order paths $(\mathfrak{R}^{(2)}_{k})_{k\in m+1}$ in $\mathrm{Pth}_{\boldsymbol{\mathcal{A}}^{(2)},s}^{m+1}$ satisfying that 
\begin{enumerate}
\item $\mathfrak{R}^{(2)}_{0}=\mathfrak{P}^{(2)}$;
\item $\mathfrak{R}^{(2)}_{m}=\mathfrak{Q}^{(2)}$;
\item for every $k\in m$, 
$(\mathfrak{R}^{(2)}_{k},\mathfrak{R}^{(2)}_{k+1})\in 
\mathrm{Ker}(\mathrm{CH}^{(2)})_{s}\cup \Upsilon^{[1]}_{s}$.
\end{enumerate}

We prove the statement by induction on $m$.

\textsf{Base step of the induction.}

For the case $m=1$, we have that either the pair $(\mathfrak{P}^{(2)},\mathfrak{Q}^{(2)})$ is in $\mathrm{Ker}(\mathrm{CH}^{(2)})_{s}$ or in $\Upsilon^{[1]}_{s}$. The statement follows by Lemmas~\ref{LDCH} and~\ref{LDUpsCong}, respectively.

\textsf{Inductive step of the induction.}

Assume that the statement holds for sequences of length $m$, i.e., if $(\mathfrak{R}^{(2)}_{k})_{k\in m+1}$ is a sequence of length $m$ in $\mathrm{Pth}_{\boldsymbol{\mathcal{A}}^{(2)},s}^{m+1}$ satisfying that 
\begin{enumerate}
\item $\mathfrak{R}^{(2)}_{0}=\mathfrak{P}^{(2)}$;
\item $\mathfrak{R}^{(2)}_{m}=\mathfrak{Q}^{(2)}$;
\item for every $k\in m$, 
$(\mathfrak{R}^{(2)}_{k},\mathfrak{R}^{(2)}_{k+1})\in 
\mathrm{Ker}(\mathrm{CH}^{(2)})_{s}\cup \Upsilon^{[1]}_{s}$,
\end{enumerate}
then it holds that
\begin{enumerate}
\item[(i)] $\bb{\mathfrak{P}^{(2)}}=\bb{\mathfrak{Q}^{(2)}}$;
\item[(ii)] $\mathrm{sc}_{s}^{([1],2)}(\mathfrak{P}^{(2)})=\mathrm{sc}_{s}^{([1],2)}(\mathfrak{Q}^{(2)})$;
\item[(iii)] $\mathrm{tg}_{s}^{([1],2)}(\mathfrak{P}^{(2)})=\mathrm{tg}_{s}^{([1],2)}(\mathfrak{Q}^{(2)})$.
\end{enumerate}

Now we prove it for a sequence of length $m+1$. Let $(\mathfrak{R}^{(2)}_{k})_{k\in m+2}$ be a sequence of length $m+1$ in $\mathrm{Pth}_{\boldsymbol{\mathcal{A}}^{(2)},s}^{m+2}$ satisfying that 
\begin{enumerate}
\item $\mathfrak{R}^{(2)}_{0}=\mathfrak{P}^{(2)}$;
\item $\mathfrak{R}^{(2)}_{m+1}=\mathfrak{Q}^{(2)}$;
\item for every $k\in m+1$, 
$(\mathfrak{R}^{(2)}_{k},\mathfrak{R}^{(2)}_{k+1})\in 
\mathrm{Ker}(\mathrm{CH}^{(2)})_{s}\cup \Upsilon^{[1]}_{s}$.
\end{enumerate}

Consider the initial subsequence $(\mathfrak{R}^{(2)}_{k})_{k\in m+1}$. It is a sequence of length $m$ in $\mathrm{Pth}_{\boldsymbol{\mathcal{A}}^{(2)},s}^{m+1}$ satisfying that 
\begin{enumerate}
\item $\mathfrak{R}^{(2)}_{0}=\mathfrak{P}^{(2)}$;
\item $\mathfrak{R}^{(2)}_{m}=\mathfrak{R}^{(2)}_{m}$;
\item for every $k\in m$, 
$(\mathfrak{R}^{(2)}_{k},\mathfrak{R}^{(2)}_{k+1})\in 
\mathrm{Ker}(\mathrm{CH}^{(2)})_{s}\cup \Upsilon^{[1]}_{s}$.
\end{enumerate}
By the inductive hypothesis it holds that 
\begin{enumerate}
\item[(i)] $\bb{\mathfrak{P}^{(2)}}=\bb{\mathfrak{R}^{(2)}_{m}}$;
\item[(ii)] $\mathrm{sc}_{s}^{([1],2)}(\mathfrak{P}^{(2)})=\mathrm{sc}_{s}^{([1],2)}(\mathfrak{R}^{(2)}_{m})$;
\item[(iii)] $\mathrm{tg}_{s}^{([1],2)}(\mathfrak{P}^{(2)})=\mathrm{tg}_{s}^{([1],2)}(\mathfrak{R}^{(2)}_{m})$.
\end{enumerate}

Now consider the final subsequence $(\mathfrak{R}^{(2)}_{m},\mathfrak{Q}^{(2)})$. It is a sequence of length $1$ in $\mathrm{Pth}_{\boldsymbol{\mathcal{A}}^{(2)},s}^{2}$ from $\mathfrak{R}^{(2)}_{m}$ to $\mathfrak{Q}^{(2)}$ satisfying that 
$(\mathfrak{R}^{(2)}_{m},\mathfrak{Q}^{(2)})\in 
\mathrm{Ker}(\mathrm{CH}^{(2)})_{s}\cup \Upsilon^{[1]}_{s}$. By the base case, we have that 
\begin{enumerate}
\item[(i)] $\bb{\mathfrak{R}^{(2)}_{m}}=\bb{\mathfrak{Q}^{(2)}}$;
\item[(ii)] $\mathrm{sc}_{s}^{([1],2)}(\mathfrak{R}^{(2)}_{m})=\mathrm{sc}_{s}^{([1],2)}(\mathfrak{Q}^{(2)})$;
\item[(iii)] $\mathrm{tg}_{s}^{([1],2)}(\mathfrak{R}^{(2)}_{m})=\mathrm{tg}_{s}^{([1],2)}(\mathfrak{Q}^{(2)})$.
\end{enumerate}

All in all, we conclude that 
\begin{enumerate}
\item[(i)] $\bb{\mathfrak{P}^{(2)}}=\bb{\mathfrak{Q}^{(2)}}$;
\item[(ii)] $\mathrm{sc}_{s}^{([1],2)}(\mathfrak{P}^{(2)})=\mathrm{sc}_{s}^{([1],2)}(\mathfrak{Q}^{(2)})$;
\item[(iii)] $\mathrm{tg}_{s}^{([1],2)}(\mathfrak{P}^{(2)})=\mathrm{tg}_{s}^{([1],2)}(\mathfrak{Q}^{(2)})$.
\end{enumerate}

This completes the proof.
\end{proof}

\begin{restatable}{corollary}{CDV}
\label{CDV} Let $s$ be a sort in $S$ and $(\mathfrak{P}^{(2)},\mathfrak{Q}^{(2)})\in \mathrm{Ker}(\mathrm{CH}^{(2)}_{s})\vee \Upsilon^{[1]}_{s}$. Then the following statements hold
\begin{enumerate}
\item[(i)] $\mathrm{sc}^{(0,2)}(\mathfrak{P}^{(2)})=\mathrm{sc}^{(0,2)}(\mathfrak{Q}^{(2)})$;
\item[(ii)] $\mathrm{tg}^{(0,2)}(\mathfrak{P}^{(2)})=\mathrm{tg}^{(0,2)}(\mathfrak{Q}^{(2)})$.
\end{enumerate}
\end{restatable}

\begin{proof}
Unravelling the mappings $\mathrm{sc}^{(0,2)}$ and $\mathrm{tg}^{(0,2)}$, introduced in Definition~\ref{DDScTgZ}, we have that 
\begin{align*}
\mathrm{sc}^{(0,2)}&=
\mathrm{sc}^{(0,[1])}
\circ
\mathrm{ip}^{([1],X)@}
\circ
\mathrm{sc}^{([1],2)};
&
\mathrm{tg}^{(0,2)}&=
\mathrm{tg}^{(0,[1])}
\circ
\mathrm{ip}^{([1],X)@}
\circ
\mathrm{tg}^{([1],2)}.
\end{align*}

Since $(\mathfrak{P}^{(2)},\mathfrak{Q}^{(2)})\in\mathrm{Ker}(\mathrm{CH}^{(2)}_{s})\vee \Upsilon^{[1]}_{s}$, then, from Lemma~\ref{LDV}, we have that 
\begin{align*}
\mathrm{sc}^{([1],2)}_{s}\left(
\mathfrak{P}^{(2)}
\right)&
=\mathrm{sc}^{([1],2)}_{s}\left(
\mathfrak{Q}^{(2)}
\right);
&
\mathrm{tg}^{([1],2)}_{s}\left(
\mathfrak{P}^{(2)}
\right)
&=\mathrm{tg}^{([1],2)}_{s}\left(
\mathfrak{Q}^{(2)}
\right).
\end{align*} 

The statement follows directly.
\end{proof}

In the following proposition we prove that each $(2,[1])$-identity second-order path is only related to itself for the supremum relation.

\begin{restatable}{proposition}{PDVDUId}
\label{PDVDUId} Let $s$ be a sort in $S$ and $( \mathfrak{P}^{(2)} ,\mathfrak{Q}^{(2)})\in \mathrm{Ker}(\mathrm{CH}^{(2)})_{s}\vee \Upsilon^{[1]}_{s}$. If $\mathfrak{P}^{(2)}$ or $\mathfrak{Q}^{(2)}$ is a $(2,[1])$-identity second-order path, then $\mathfrak{P}^{(2)}=\mathfrak{Q}^{(2)}$.
\end{restatable}
\begin{proof}
Let us recall that, for a sort $s\in S$, two second-order paths $\mathfrak{P}^{(2)}$ and $\mathfrak{Q}^{(2)}$ in $\mathrm{Pth}_{\boldsymbol{\mathcal{A}}^{(2)},s}$ are related under $\mathrm{Ker}(\mathrm{CH}^{(2)}_{s})\vee \Upsilon^{[1]}_{s}$ if, and only if, there exists a natural number $m\in \mathbb{N}-\{0\}$ and a sequence of second-order paths $(\mathfrak{R}^{(2)}_{k})_{k\in m+1}$ in $\mathrm{Pth}_{\boldsymbol{\mathcal{A}}^{(2)},s}^{m+1}$ satisfying that 
\begin{enumerate}
\item $\mathfrak{R}^{(2)}_{0}=\mathfrak{P}^{(2)}$;
\item $\mathfrak{R}^{(2)}_{m}=\mathfrak{Q}^{(2)}$;
\item for every $k\in m$, 
$(\mathfrak{R}^{(2)}_{k},\mathfrak{R}^{(2)}_{k+1})\in 
\mathrm{Ker}(\mathrm{CH}^{(2)})_{s}\cup \Upsilon^{[1]}_{s}$.
\end{enumerate}

We prove the statement by induction on $m$.

\textsf{Base step of the induction.}

For the case $m=1$, we have that either the pair $(\mathfrak{P}^{(2)},\mathfrak{Q}^{(2)})$ is in $\mathrm{Ker}(\mathrm{CH}^{(2)})_{s}$ or in $\Upsilon^{[1]}_{s}$. The statement follows by Corollaries~\ref{CDCHUId} and~\ref{CDUpsCongDUId}, respectively.

\textsf{Inductive step of the induction.}

Assume that the statement holds for sequences of length $m$, i.e., if $(\mathfrak{R}^{(2)}_{k})_{k\in m+1}$ is a sequence of length $m$ in $\mathrm{Pth}_{\boldsymbol{\mathcal{A}}^{(2)},s}^{m+1}$ satisfying that 
\begin{enumerate}
\item $\mathfrak{R}^{(2)}_{0}=\mathfrak{P}^{(2)}$;
\item $\mathfrak{R}^{(2)}_{m}=\mathfrak{Q}^{(2)}$;
\item for every $k\in m$, 
$(\mathfrak{R}^{(2)}_{k},\mathfrak{R}^{(2)}_{k+1})\in 
\mathrm{Ker}(\mathrm{CH}^{(2)})_{s}\cup \Upsilon^{[1]}_{s}$,
\end{enumerate}
then, if $\mathfrak{P}^{(2)}$ or $\mathfrak{Q}^{(2)}$ is a $(2,[1])$-identity second-order path, it holds that  $\mathfrak{P}^{(2)}=\mathfrak{Q}^{(2)}$.

Now we prove it for a sequence of length $m+1$. Let $(\mathfrak{R}^{(2)}_{k})_{k\in m+2}$ be a sequence of length $m+1$ in $\mathrm{Pth}_{\boldsymbol{\mathcal{A}}^{(2)},s}^{m+2}$ satisfying that 
\begin{enumerate}
\item $\mathfrak{R}^{(2)}_{0}=\mathfrak{P}^{(2)}$;
\item $\mathfrak{R}^{(2)}_{m+1}=\mathfrak{Q}^{(2)}$;
\item for every $k\in m+1$, 
$(\mathfrak{R}^{(2)}_{k},\mathfrak{R}^{(2)}_{k+1})\in 
\mathrm{Ker}(\mathrm{CH}^{(2)})_{s}\cup \Upsilon^{[1]}_{s}$.
\end{enumerate}

Assume without loss of generality that $\mathfrak{P}^{(2)}$ is a $(2,[1])$-identity second-order path. Consider the initial subsequence $(\mathfrak{R}^{(2)}_{k})_{k\in m+1}$. It is a sequence of length $m$ in $\mathrm{Pth}_{\boldsymbol{\mathcal{A}}^{(2)},s}^{m+1}$ satisfying that 
\begin{enumerate}
\item $\mathfrak{R}^{(2)}_{0}=\mathfrak{P}^{(2)}$;
\item $\mathfrak{R}^{(2)}_{m}=\mathfrak{R}^{(2)}_{m}$;
\item for every $k\in m$, 
$(\mathfrak{R}^{(2)}_{k},\mathfrak{R}^{(2)}_{k+1})\in 
\mathrm{Ker}(\mathrm{CH}^{(2)})_{s}\cup \Upsilon^{[1]}_{s}$.
\end{enumerate}
By the inductive hypothesis it holds that $\mathfrak{P}^{(2)}=\mathfrak{R}^{(2)}_{m}$.

Now consider the final subsequence $(\mathfrak{R}^{(2)}_{m},\mathfrak{Q}^{(2)})$. It is a sequence of length $1$ in $\mathrm{Pth}_{\boldsymbol{\mathcal{A}}^{(2)},s}^{2}$ from $\mathfrak{R}^{(2)}_{m}$ to $\mathfrak{Q}^{(2)}$ satisfying that 
$(\mathfrak{R}^{(2)}_{m},\mathfrak{Q}^{(2)})\in 
\mathrm{Ker}(\mathrm{CH}^{(2)})_{s}\cup \Upsilon^{[1]}_{s}$. By the base case, we have that $\mathfrak{R}^{(2)}_{m}=\mathfrak{Q}^{(2)}$.

All in all, we conclude that  $\mathfrak{P}^{(2)}=\mathfrak{Q}^{(2)}$.

This completes the proof.
\end{proof}

We now introduce different many-sorted mappings related with the $S$-sorted quotient $\llbracket\mathrm{Pth}_{\boldsymbol{\mathcal{A}}^{(2)}} \rrbracket$.

\begin{restatable}{definition}{DDVQuot}
\label{DDVQuot} We define the projection, denoted by $\mathrm{pr}^{\llbracket \cdot \rrbracket}$, to be the $S$-sorted mapping 
\[
\begin{array}{rccl}
\mathrm{pr}^{\llbracket \cdot\rrbracket }\colon&\mathrm{Pth}_{\boldsymbol{\mathcal{A}}^{(2)}}
&\longrightarrow&\llbracket \mathrm{Pth}_{\boldsymbol{\mathcal{A}}^{(2)}}
 \rrbracket 
\end{array}
\]
\index{projection!second-order!$\mathrm{pr}^{\llbracket \cdot\rrbracket }$}
that, for every sort $s\in S$, maps a second-order path $\mathfrak{P}^{(2)}$ in $\mathrm{Pth}_{\boldsymbol{\mathcal{A}}^{(2)},s}$ to $\llbracket\mathfrak{P}^{(2)}\rrbracket_{s}$, its equivalence class under the supremum of the kernel of the second-order Curry-Howard mapping and the  congruence $\Upsilon^{[1]}$.
\end{restatable}

We next introduce the notation for the composition of the embeddings to $\mathrm{Pth}_{\boldsymbol{\mathcal{A}}^{(2)}}$ with the projection to the supremum of the kernel of the second-order Curry-Howard mapping and the congruence $\Upsilon^{[1]}$.

\begin{restatable}{definition}{DDVEch}
\label{DDVEch} We will denote by
\begin{enumerate}
\item $\mathrm{ip}^{(\llbracket 2\rrbracket,X)}$ the $S$-sorted mapping from $X$ to $\llbracket \mathrm{Pth}_{\boldsymbol{\mathcal{A}}^{(2)}}\rrbracket$ given by the composition $\mathrm{ip}^{(\llbracket 2\rrbracket,X)}=\mathrm{pr}^{\llbracket \cdot \rrbracket}\circ\mathrm{ip}^{(2,X)}$, i.e., for every sort $s\in S$, $\mathrm{ip}^{(\llbracket 2\rrbracket,X)}_{s}$ sends a variable $x\in X_{s}$ to the  class $\llbracket \mathrm{ip}^{(2,X)}_{s}(x)\rrbracket_{s}$ in $\llbracket \mathrm{Pth}_{\boldsymbol{\mathcal{A}}^{(2)}} \rrbracket_{s}$.
\index{identity!second-order!$\mathrm{ip}^{(\llbracket 2\rrbracket,X)}$}
\item $\mathrm{ech}^{(\llbracket 2\rrbracket,\mathcal{A})}$ the $S$-sorted mapping from $\mathcal{A}$ to $\llbracket \mathrm{Pth}_{\boldsymbol{\mathcal{A}}^{(2)}}\rrbracket$ given by the composition $\mathrm{ech}^{(\llbracket 2\rrbracket,\mathcal{A})}=\mathrm{pr}^{\llbracket \cdot \rrbracket}\circ\mathrm{ech}^{(2,\mathcal{A})}$, i.e., for every sort $s\in S$, $\mathrm{ech}^{(\llbracket 2\rrbracket,\mathcal{A})}_{s}$ sends a rewrite rule $\mathfrak{p}\in \mathcal{A}_{s}$ to the  class $\llbracket \mathrm{ech}^{(2,\mathcal{A})}_{s}(\mathfrak{p})\rrbracket_{s}$ in $\llbracket \mathrm{Pth}_{\boldsymbol{\mathcal{A}}^{(2)}} \rrbracket_{s}$
\index{echelon!second-order!$\mathrm{ech}^{(\llbracket 2\rrbracket,\mathcal{A})}$}
\item $\mathrm{ech}^{(\llbracket 2\rrbracket,\mathcal{A}^{(2)})}$ the $S$-sorted mapping from $\mathcal{A}^{(2)}$ to $\llbracket \mathrm{Pth}_{\boldsymbol{\mathcal{A}}^{(2)}}\rrbracket$ given by the composition $\mathrm{ech}^{(\llbracket 2\rrbracket,\mathcal{A}^{(2)})}=\mathrm{pr}^{\llbracket \cdot \rrbracket}\circ\mathrm{ech}^{(2,\mathcal{A}^{(2)})}$, i.e., for every sort $s\in S$, $\mathrm{ech}^{(\llbracket 2\rrbracket,\mathcal{A}^{(2)})}_{s}$ sends a second-order rewrite rule $\mathfrak{p}^{(2)}\in \mathcal{A}^{(2)}_{s}$ to the  class $\llbracket \mathrm{ech}^{(2,\mathcal{A}^{(2)})}_{s}(\mathfrak{p}^{(2)})\rrbracket_{s}$ in $\llbracket \mathrm{Pth}_{\boldsymbol{\mathcal{A}}^{(2)}} \rrbracket_{s}$.
\index{echelon!second-order!$\mathrm{ech}^{(\llbracket 2\rrbracket,\mathcal{A}^{(2)})}$}
\end{enumerate}
The above $S$-sorted mappings are depicted in the diagram of Figure~\ref{FDVEch}.
\end{restatable}

\begin{figure}
\begin{tikzpicture}
[ACliment/.style={-{To [angle'=45, length=5.75pt, width=4pt, round]}},scale=.8]
\node[] (x) at (0,0) [] {$X$};
\node[] (a) at (0,-1.5) [] {$\mathcal{A}$};
\node[] (a2) at (0,-3) [] {$\mathcal{A}^{(2)}$};
\node[] (T) at (6,-3) [] {$
\llbracket \mathrm{Pth}_{\boldsymbol{\mathcal{A}}^{(2)}}\rrbracket$};
\draw[ACliment, bend left=20]  (x) to node [above right] {$\mathrm{ip}^{(\llbracket 2\rrbracket,X)}$} (T);
\draw[ACliment, bend left=10]  (a) to node [midway, fill=white] {$\mathrm{ech}^{(\llbracket 2\rrbracket,\mathcal{A})}$} (T);
\draw[ACliment]  (a2) to node [below] {$\mathrm{ech}^{(\llbracket 2\rrbracket,\mathcal{A}^{(2)})}$} (T);
\end{tikzpicture}
\caption{Quotient path embeddings relative to $X$, $\mathcal{A}$ and $\mathcal{A}^{(2)}$.}\label{FDVEch}
\end{figure}

Finally, we introduce some terminology for the $S$-sorted mappings of sources, targets and identity second-order paths at each layer with with respect to the many-sorted first-order quotient set of second-order paths with respect to the supremum of the kernel of the second-order Curry-Howard mapping and the congruence $\Upsilon^{[1]}$.

\begin{restatable}{definition}{DDVDU}
\label{DDVDU} We will denote by
\begin{enumerate}
\item  $\mathrm{sc}^{([1],\llbracket 2\rrbracket)}$ the $S$-sorted mapping from $\llbracket \mathrm{Pth}_{\boldsymbol{\mathcal{A}}^{(2)}}\rrbracket$ to $[\mathrm{PT}_{\boldsymbol{\mathcal{A}}}]$ that, for every $s\in S$, assigns to an equivalence class $\llbracket \mathfrak{P}^{(2)}\rrbracket_{s}$ in $\llbracket \mathrm{Pth}_{\boldsymbol{\mathcal{A}}^{(2)}} \rrbracket_{s}$ with $\mathfrak{P}^{(2)}\in\mathrm{Pth}_{\boldsymbol{\mathcal{A}}^{(2)},s}$ the path term class $\mathrm{sc}^{([1],2)}_{s}(\mathfrak{P}^{(2)})$, i.e., the value of the mapping $\mathrm{sc}^{([1],2)}_{s}$ at any   equivalence class representative. That $\mathrm{sc}^{([1],\llbracket 2\rrbracket)}$ is well-defined follows from Lemma~\ref{LDV};
\index{source!second-order!$\mathrm{sc}^{([1],\llbracket 2\rrbracket)}$}
\item $\mathrm{tg}^{([1],\llbracket 2\rrbracket)}$ the $S$-sorted mapping  from $\llbracket \mathrm{Pth}_{\boldsymbol{\mathcal{A}}^{(2)}}\rrbracket$ to $[\mathrm{PT}_{\boldsymbol{\mathcal{A}}}]$ that, for every $s\in S$, assigns to an equivalence class $\llbracket \mathfrak{P}^{(2)}\rrbracket_{s}$ in $\llbracket \mathrm{Pth}_{\boldsymbol{\mathcal{A}}^{(2)}} \rrbracket_{s}$ with $\mathfrak{P}^{(2)}\in\mathrm{Pth}_{\boldsymbol{\mathcal{A}}^{(2)},s}$ the path term class $\mathrm{tg}^{([1],2)}_{s}(\mathfrak{P}^{(2)})$, i.e., the value of the mapping $\mathrm{tg}^{([1],2)}_{s}$ at any  equivalence class representative. That $\mathrm{tg}^{([1],\llbracket 2\rrbracket)}$ is well-defined follows from Lemma~\ref{LDV};
\index{target!second-order!$\mathrm{tg}^{([1],\llbracket 2\rrbracket)}$}
\item $\mathrm{ip}^{(\llbracket 2\rrbracket,[1])\sharp}$  the $S$-sorted mapping from $[\mathrm{PT}_{\boldsymbol{\mathcal{A}}}]$ to $\llbracket \mathrm{Pth}_{\boldsymbol{\mathcal{A}}^{(2)}}\rrbracket$ given by the composition
$\mathrm{ip}^{(\llbracket 2\rrbracket,[1])\sharp}=\mathrm{pr}^{\llbracket \cdot \rrbracket }\circ\mathrm{ip}^{(2,[1])\sharp}$. That is, for every sort $s\in S$, $\mathrm{ip}^{(\llbracket 2\rrbracket,[1])\sharp}$ assigns to a path term  class $[P]_{s}$ in $[\mathrm{PT}_{\boldsymbol{\mathcal{A}}}]_{s}$ the class with respect to supremum of the second-order Curry-Howard kernel and the congruence $\Upsilon^{[1]}$ of the $(2,[1])$-identity second-order path $\mathrm{ip}^{(2,[1])\sharp}_{s}([P]_{s})$, i.e., 
$
\mathrm{ip}^{(\llbracket 2\rrbracket,[1])\sharp}_{s}(
[P]_{s}
)=\llbracket \mathrm{ip}^{(2,[1])\sharp}_{s}([P]_{s})\rrbracket_{s}.
$
\index{identity!second-order!$\mathrm{ip}^{(\llbracket 2\rrbracket,[1])\sharp}$}
\end{enumerate}
\end{restatable}

\begin{proposition} The following equalities hold
\begin{itemize}
\item[(i)] $\mathrm{sc}^{([1],\llbracket 2\rrbracket)}\circ \mathrm{ip}^{(\llbracket 2\rrbracket,[1])\sharp}=\mathrm{id}^{
[
\mathrm{PT}_{\boldsymbol{\mathcal{A}}}
]
}$;
\item[(ii)] $\mathrm{tg}^{([1],\llbracket 2\rrbracket)}\circ \mathrm{ip}^{(\llbracket 2\rrbracket,[1])\sharp}=\mathrm{id}^{
[
\mathrm{PT}_{\boldsymbol{\mathcal{A}}}
]
}$.
\end{itemize}
\end{proposition}

\begin{proposition} The following equalities hold
\begin{itemize}
\item[(i)] $\mathrm{sc}^{([1],\llbracket 2\rrbracket)}\circ \mathrm{ip}^{(\llbracket 2\rrbracket,X)}
=\eta^{([1],X)};$
\item[(ii)] $\mathrm{tg}^{([1],\llbracket 2\rrbracket)}\circ \mathrm{ip}^{(\llbracket 2\rrbracket,X)}
=\eta^{([1],X)};$
\item[(iii)] $\mathrm{ip}^{(\llbracket 2\rrbracket,[1])\sharp}\circ \eta^{([1],X)}=
\mathrm{ip}^{(\llbracket 2\rrbracket,X)}.
$
\end{itemize}

The reader is advised to consult the diagram appearing in Figure~\ref{FDVDU}.
\end{proposition}

\begin{figure}
\begin{center}
\begin{tikzpicture}
[ACliment/.style={-{To [angle'=45, length=5.75pt, width=4pt, round]}},scale=1]
\node[] (x) at (0,0) [] {$X$};
\node[] (pt) at (6,0) [] {$[\mathrm{PT}_{\boldsymbol{\mathcal{A}}}]$};
\node[] (2pq) at (6,-3) []  {$\llbracket \mathrm{Pth}_{\boldsymbol{\mathcal{A}}^{(2)}} \rrbracket 
$};

\draw[ACliment]  (x) 	to node [above right]	
{$\eta^{([1],X)}$} (pt);
\draw[ACliment, bend right=10]  (x) 	to node [below left]	
{$\mathrm{ip}^{(\llbracket 2\rrbracket,X)}$} (2pq);

\node[] (B1) at (6,-1.5)  [] {};
\draw[ACliment]  ($(B1)+(0,1.2)$) to node [above, fill=white] {
$\mathrm{ip}^{(\llbracket 2\rrbracket,[1])\sharp}$
} ($(B1)+(0,-1.2)$);
\draw[ACliment, bend right]  ($(B1)+(.5,-1.2)$) to node [ below, fill=white] {
$ \mathrm{tg}^{([1],\llbracket 2\rrbracket)}$
} ($(B1)+(.5,1.2)$);
\draw[ACliment, bend left]  ($(B1)+(-.5,-1.2)$) to node [below, fill=white] {
$ \mathrm{sc}^{([1],\llbracket 2\rrbracket)}$
} ($(B1)+(-.5,1.2)$);

\end{tikzpicture}
\end{center}
\caption{Quotient mappings relative to $X$ at layers 1 \& 2.}
\label{FDVDU}
\end{figure}

\begin{proposition} The following equalities hold
\begin{itemize}
\item[(i)] $\mathrm{sc}^{([1],\llbracket 2\rrbracket)}\circ \mathrm{ech}^{(\llbracket 2\rrbracket,\mathcal{A})}
=\eta^{([1],\mathcal{A})};$
\item[(ii)] $\mathrm{tg}^{([1],\llbracket 2\rrbracket)}\circ \mathrm{ech}^{(\llbracket 2\rrbracket,\mathcal{A})}
= \eta^{([1],\mathcal{A})};$
\item[(iii)] $\mathrm{ip}^{(\llbracket 2\rrbracket,[1])\sharp}\circ \eta^{([1],\mathcal{A})}=
\mathrm{ech}^{(\llbracket 2\rrbracket,\mathcal{A})}.
$
\end{itemize}

The reader is advised to consult the diagram appearing in Figure~\ref{FDVDA}.
\end{proposition}

\begin{figure}
\begin{center}
\begin{tikzpicture}
[ACliment/.style={-{To [angle'=45, length=5.75pt, width=4pt, round]}},scale=1]
\node[] (x) at (0,0) [] {$\mathcal{A}$};
\node[] (pt) at (6,0) [] {$[\mathrm{PT}_{\boldsymbol{\mathcal{A}}}]$};
\node[] (2pq) at (6,-3) []  {$\llbracket \mathrm{Pth}_{\boldsymbol{\mathcal{A}}^{(2)}}
\rrbracket 
$};

\draw[ACliment]  (x) 	to node [above right]	
{$\eta^{([1],\mathcal{A})}$} (pt);
\draw[ACliment, bend right=10]  (x) 	to node [below left]	
{$\mathrm{ech}^{(\llbracket 2\rrbracket,\mathcal{A})}$} (2pq);

\node[] (B1) at (6,-1.5)  [] {};
\draw[ACliment]  ($(B1)+(0,1.2)$) to node [above, fill=white] {
$\mathrm{ip}^{(\llbracket 2\rrbracket,[1])\sharp}$
} ($(B1)+(0,-1.2)$);
\draw[ACliment, bend right]  ($(B1)+(.5,-1.2)$) to node [ below, fill=white] {
$ \mathrm{tg}^{([1],\llbracket 2\rrbracket)}$
} ($(B1)+(.5,1.2)$);
\draw[ACliment, bend left]  ($(B1)+(-.5,-1.2)$) to node [below, fill=white] {
$ \mathrm{sc}^{([1],\llbracket 2\rrbracket)}$
} ($(B1)+(-.5,1.2)$);

\end{tikzpicture}
\end{center}
\caption{Many-sorted quotient mappings relative to $\mathcal{A}$ at layers 1 \& 2.}
\label{FDVDA}
\end{figure}

\begin{restatable}{definition}{DDVDZ}
\label{DDVDZ} We will denote by
\begin{enumerate}
\item $\mathrm{sc}^{(0,\llbracket 2\rrbracket)}$ the $S$-sorted mapping from $\llbracket \mathrm{Pth}_{\boldsymbol{\mathcal{A}}^{(2)}}\rrbracket$ to $\mathrm{T}_{\Sigma}(X)$ that, for every $s\in S$, assigns to an equivalence class $\llbracket \mathfrak{P}^{(2)}\rrbracket_{s}$ in $\llbracket \mathrm{Pth}_{\boldsymbol{\mathcal{A}}^{(2)}} \rrbracket_{s}$ with $\mathfrak{P}^{(2)}\in\mathrm{Pth}_{\boldsymbol{\mathcal{A}}^{(2)},s}$ the  term $\mathrm{sc}^{(0,2)}_{s}(\mathfrak{P}^{(2)})$, i.e., the value of the mapping $\mathrm{sc}^{(0,2)}_{s}$ at any equivalence class representative. That $\mathrm{sc}^{(0,2)}$ is well-defined follows from Corollary~\ref{CDV};
\index{source!second-order!$\mathrm{sc}^{(0,\llbracket 2\rrbracket)}$}

\item $\mathrm{tg}^{(0,\llbracket 2\rrbracket)}$ the $S$-sorted mapping from $\llbracket \mathrm{Pth}_{\boldsymbol{\mathcal{A}}^{(2)}}\rrbracket$ to $\mathrm{T}_{\Sigma}(X)$ that, for every $s\in S$, assigns to an equivalence class $\llbracket \mathfrak{P}^{(2)}\rrbracket_{s}$ in $\llbracket \mathrm{Pth}_{\boldsymbol{\mathcal{A}}^{(2)}} \rrbracket_{s}$ with $\mathfrak{P}^{(2)}\in\mathrm{Pth}_{\boldsymbol{\mathcal{A}}^{(2)},s}$ the  term $\mathrm{tg}^{(0,2)}_{s}(\mathfrak{P}^{(2)})$, i.e., the value of the mapping $\mathrm{tg}^{(0,2)}_{s}$ at any equivalence class representative. That $\mathrm{tg}^{(0,2)}$ is well-defined follows from Corollary~\ref{CDV};
\index{target!second-order!$\mathrm{tg}^{(0,\llbracket 2\rrbracket)}$}
\item $\mathrm{ip}^{(\llbracket 2\rrbracket,0)\sharp}$  the $S$-sorted mapping from $\mathrm{T}_{\Sigma}(X)$ to $\llbracket \mathrm{Pth}_{\boldsymbol{\mathcal{A}}^{(2)}}\rrbracket$ given by the composition
$\mathrm{ip}^{(\llbracket 2\rrbracket,0)\sharp}=\mathrm{pr}^{\llbracket \cdot \rrbracket}\circ\mathrm{ip}^{(2,0)\sharp}$. That is, for every sort $s\in S$, $\mathrm{ip}^{(\llbracket 2\rrbracket,0)\sharp}$ assigns to a term $P$ in $\mathrm{T}_{\Sigma}(X)_{s}$ the class with respect to the supremum of the kernel of second-order Curry-Howard mapping and the  congruence $\Upsilon^{[1]}$ of the $(2,0)$-identity second-order path $\mathrm{ip}^{(2,0)\sharp}_{s}(P)$, i.e., 
$
\mathrm{ip}^{(\llbracket 2\rrbracket,0)\sharp}_{s}(
P
)=\llbracket \mathrm{ip}^{(2,0)\sharp}_{s}(P)\rrbracket _{s}.
$
\index{identity!second-order!$\mathrm{ip}^{(\llbracket 2\rrbracket,0)\sharp}$}
\end{enumerate} 
\end{restatable}

\begin{proposition} The following equalities hold
\begin{itemize}
\item[(i)] $\mathrm{sc}^{(0,\llbracket 2\rrbracket)}=\mathrm{sc}^{(0,[1])}\circ \mathrm{sc}^{([1],\llbracket 2\rrbracket)}$;
\item[(ii)] $\mathrm{tg}^{(0,\llbracket 2\rrbracket)}=\mathrm{tg}^{(0,[1])}\circ\mathrm{tg}^{([1],\llbracket 2\rrbracket)}$;
\item[(iii)] $\mathrm{ip}^{(\llbracket 2\rrbracket,0)\sharp}=\mathrm{ip}^{(\llbracket 2\rrbracket,[1])\sharp}\circ\mathrm{ip}^{([1],0)\sharp}$.
\end{itemize}
\end{proposition}

\begin{proposition} The following equalities hold
\begin{itemize}
\item[(i)] $\mathrm{sc}^{(0,\llbracket 2\rrbracket)}\circ \mathrm{ip}^{(\llbracket 2\rrbracket,0)\sharp}=\mathrm{id}^{
\mathrm{T}_{\Sigma}(X)
}$;
\item[(ii)] $\mathrm{tg}^{(0,\llbracket 2\rrbracket)}\circ\mathrm{ip}^{(\llbracket 2\rrbracket,0)\sharp}=\mathrm{id}^{
\mathrm{T}_{\Sigma}(X)
}$.
\end{itemize}
\end{proposition}

\begin{proposition} The following equalities hold
\begin{itemize}
\item[(i)] $\mathrm{sc}^{(0,\llbracket 2\rrbracket)}\circ \mathrm{ip}^{(\llbracket 2\rrbracket,X)}
=\eta^{(0,X)};$
\item[(ii)] $\mathrm{tg}^{(0,\llbracket 2\rrbracket)}\circ \mathrm{ip}^{(\llbracket 2\rrbracket,X)}
=\eta^{(0,X)};$
\item[(iii)] $\mathrm{ip}^{(\llbracket 2\rrbracket,0)\sharp}\circ\eta^{(0,X)}=
\mathrm{ip}^{(\llbracket 2\rrbracket,X)}.
$
\end{itemize}

The reader is advised to consult the diagram appearing in Figure~\ref{FDVDZ}.
\end{proposition}

\begin{figure}
\begin{center}
\begin{tikzpicture}
[ACliment/.style={-{To [angle'=45, length=5.75pt, width=4pt, round]}},scale=1]
\node[] (x) at (0,0) [] {$X$};
\node[] (pt) at (6,0) [] {$\mathrm{T}_{\Sigma}(X)$};
\node[] (2pq) at (6,-3) []  {$\llbracket \mathrm{Pth}_{\boldsymbol{\mathcal{A}}^{(2)}}\rrbracket 
$};

\draw[ACliment]  (x) 	to node [above right]	
{$\eta^{(0,X)}$} (pt);
\draw[ACliment, bend right=10]  (x) 	to node [below left]	
{$\mathrm{ip}^{(\llbracket 2\rrbracket,X)}$} (2pq);

\node[] (B1) at (6,-1.5)  [] {};
\draw[ACliment]  ($(B1)+(0,1.2)$) to node [above, fill=white] {
$\mathrm{ip}^{(\llbracket 2\rrbracket,0)\sharp}$
} ($(B1)+(0,-1.2)$);
\draw[ACliment, bend right]  ($(B1)+(.5,-1.2)$) to node [ below, fill=white] {
$ \mathrm{tg}^{(0,\llbracket 2\rrbracket)}$
} ($(B1)+(.5,1.2)$);
\draw[ACliment, bend left]  ($(B1)+(-.5,-1.2)$) to node [below, fill=white] {
$ \mathrm{sc}^{(0,\llbracket 2\rrbracket)}$
} ($(B1)+(-.5,1.2)$);

\end{tikzpicture}
\end{center}
\caption{Quotient mappings relative to $X$ at layers 0 \& 2.}
\label{FDVDZ}
\end{figure}

\section{
\texorpdfstring
{A structure of partial $\Sigma^{\boldsymbol{\mathcal{A}}^{(2)}}$-algebra on $\llbracket \mathrm{Pth}_{\boldsymbol{\mathcal{A}}^{(2)}}\rrbracket$}
{A structure of partial algebra on the second-order quotient}
} 

In this section we introduce the algebraic structure in the quotient of the second-order paths by the supremum of the two congruences $\mathrm{Ker}(\mathrm{CH}^{(2)})$ and $\Upsilon^{[1]}$ on it.

\begin{restatable}{proposition}{PDVCong}
\label{PDVCong} $\mathrm{Ker}(\mathrm{CH}^{(2)})\vee \Upsilon^{[1]}$ is a closed $\Sigma^{\boldsymbol{\mathcal{A}}^{(2)}}$-congruence on $\mathbf{Pth}_{\boldsymbol{\mathcal{A}}^{(2)}}$.
\end{restatable}
\begin{proof}
It follows from the fact that $\mathrm{Ker}(\mathrm{CH}^{(2)})$ and $\Upsilon^{[1]}$ are closed $\Sigma^{\boldsymbol{\mathcal{A}}^{(2)}}$-congruences  on $\mathbf{Pth}_{\boldsymbol{\mathcal{A}}^{(2)}}$, according to Proposition~\ref{PDCHCong} and Definition~\ref{DDUpsCong}, respectively.
\end{proof}

Since $\mathrm{Ker}(\mathrm{CH}^{(2)})\vee \Upsilon^{[1]}$ is a $\Sigma^{\boldsymbol{\mathcal{A}}^{(2)}}$-congruence on $\mathbf{Pth}_{\boldsymbol{\mathcal{A}}^{(2)}}$ we can endow the many-sorted quotient $\llbracket \mathrm{Pth}_{\boldsymbol{\mathcal{A}}^{(2)}}\rrbracket$ with a structure of many-sorted partial $\Sigma^{\boldsymbol{\mathcal{A}}^{(2)}}$-algebra.

\begin{restatable}{proposition}{PDVDCatAlg}
\label{PDVDCatAlg} 
\index{path!second-order!$ \llbracket\mathbf{Pth}_{\boldsymbol{\mathcal{A}}^{(2)}} \rrbracket $}
The $S$-sorted set $\llbracket \mathrm{Pth}_{\boldsymbol{\mathcal{A}}^{(2)}}\rrbracket$ is equipped, in a natural way, with a structure of many-sorted partial $\Sigma^{\boldsymbol{\mathcal{A}}^{(2)}}$-algebra. We denote by $ \llbracket\mathbf{Pth}_{\boldsymbol{\mathcal{A}}^{(2)}} \rrbracket $ the corresponding $\Sigma^{\boldsymbol{\mathcal{A}}^{(2)}}$-algebra (which is a quotient of $\mathbf{Pth}_{\boldsymbol{\mathcal{A}}^{(2)}}$, the many-sorted partial $\Sigma^{\boldsymbol{\mathcal{A}}^{(2)}}$-algebra constructed in Proposition~\ref{PDPthDCatAlg}).

The mapping $\mathrm{pr}^{\llbracket \cdot \rrbracket}$, that is,
$$
\mathrm{pr}^{\llbracket \cdot \rrbracket}
\colon
\mathbf{Pth}_{\boldsymbol{\mathcal{A}}^{(2)}}
\mor
\llbracket \mathbf{Pth}_{\boldsymbol{\mathcal{A}}^{(2)}}
\rrbracket
$$
 is a closed and surjective $\Sigma^{\boldsymbol{\mathcal{A}}^{(2)}}$-homomorphism from $\mathbf{Pth}_{\boldsymbol{\mathcal{A}}^{(2)}}$ to $ \llbracket\mathbf{Pth}_{\boldsymbol{\mathcal{A}}^{(2)}} \rrbracket $.
\end{restatable}

We next consider the $\Sigma^{\boldsymbol{\mathcal{A}}}$-reduct of the many-sorted partial $\Sigma^{\boldsymbol{\mathcal{A}}^{(2)}}$-algebra of second-order path classes and study its connections with the mappings presented in the previous section.

\begin{definition} For the partial $\Sigma^{\boldsymbol{\mathcal{A}}^{(2)}}$-algebra $ \llbracket\mathbf{Pth}_{\boldsymbol{\mathcal{A}}^{(2)}} \rrbracket $, we denote by $\llbracket \mathbf{Pth}^{(1,2)}_{\boldsymbol{\mathcal{A}}^{(2)}}
\rrbracket$ the $\Sigma^{\boldsymbol{\mathcal{A}}}$-algebra $\mathbf{in}^{\Sigma,(1,2)}\left(
\llbracket \mathbf{Pth}_{\boldsymbol{\mathcal{A}}^{(2)}}
\rrbracket\right)$. We will call  $\llbracket \mathbf{Pth}^{(1,2)}_{\boldsymbol{\mathcal{A}}^{(2)}}
\rrbracket$ the $\Sigma^{\boldsymbol{\mathcal{A}}}$-reduct of the partial $\Sigma^{\boldsymbol{\mathcal{A}}^{(2)}}$-algebra  $ \llbracket\mathbf{Pth}_{\boldsymbol{\mathcal{A}}^{(2)}} \rrbracket $.
\end{definition}

\begin{proposition} The mapping $\mathrm{pr}^{\llbracket \cdot \rrbracket}$
is a closed and surjective $\Sigma^{\boldsymbol{\mathcal{A}}}$-homomorphism from $\mathbf{Pth}^{(1,2)}_{\boldsymbol{\mathcal{A}}^{(2)}}$ to $\llbracket \mathbf{Pth}^{(1,2)}_{\boldsymbol{\mathcal{A}}^{(2)}}
\rrbracket$.
\end{proposition}

\begin{restatable}{proposition}{PDVDU}
\label{PDVDU} The many-sorted mappings $\mathrm{sc}^{([1],\llbracket 2\rrbracket)}$ and $\mathrm{tg}^{([1],\llbracket 2\rrbracket)}$ are $\Sigma^{\boldsymbol{\mathcal{A}}}$-homomorphisms from $\llbracket \mathbf{Pth}^{(1,2)}_{\boldsymbol{\mathcal{A}}^{(2)}}
\rrbracket$ to $[\mathbf{PT}_{\boldsymbol{\mathcal{A}}}]$.
\end{restatable}

\begin{restatable}{proposition}{PDVDUIp}
\label{PDVDUIp} The many-sorted mapping $\mathrm{ip}^{(\llbracket 2\rrbracket,[1])\sharp}$  is a $\Sigma^{\boldsymbol{\mathcal{A}}}$-homomorphism from  $[\mathbf{PT}_{\boldsymbol{\mathcal{A}}}]$ to $\llbracket \mathbf{Pth}^{(1,2)}_{\boldsymbol{\mathcal{A}}^{(2)}}
\rrbracket$.
\end{restatable}

We next consider the $\Sigma$-reduct of the many-sorted partial $\Sigma^{\boldsymbol{\mathcal{A}}^{(2)}}$-algebra of second-order path classes and study its connections with the mappings presented in the previous section.

\begin{definition} For the partial $\Sigma^{\boldsymbol{\mathcal{A}}^{(2)}}$-algebra $ \llbracket\mathbf{Pth}_{\boldsymbol{\mathcal{A}}^{(2)}} \rrbracket $, we denote by $\llbracket \mathbf{Pth}^{(0,2)}_{\boldsymbol{\mathcal{A}}^{(2)}}
\rrbracket$ the $\Sigma$-algebra $\mathbf{in}^{\Sigma,(0,2)}\left(
\llbracket \mathbf{Pth}_{\boldsymbol{\mathcal{A}}^{(2)}}
\rrbracket\right).$ We will call  $\llbracket \mathbf{Pth}^{(0,2)}_{\boldsymbol{\mathcal{A}}^{(2)}}
\rrbracket$ the $\Sigma$-reduct of the partial $\Sigma^{\boldsymbol{\mathcal{A}}^{(2)}}$-algebra  $ \llbracket\mathbf{Pth}_{\boldsymbol{\mathcal{A}}^{(2)}} \rrbracket $.
\end{definition}

\begin{proposition}\label{PDVDZPr} The mapping $\mathrm{pr}^{\llbracket \cdot \rrbracket}$
is a surjective $\Sigma$-homomorphism from $\mathbf{Pth}^{(0,2)}_{\boldsymbol{\mathcal{A}}^{(2)}}$ to $\llbracket \mathbf{Pth}^{(0,2)}_{\boldsymbol{\mathcal{A}}^{(2)}}
\rrbracket$.
\end{proposition}

\begin{restatable}{proposition}{PDVDZ}
\label{PDVDZ} The many-sorted mappings $\mathrm{sc}^{(0,\llbracket 2\rrbracket)}$ and $\mathrm{tg}^{(0,\llbracket 2\rrbracket)}$ are $\Sigma$-homomorphisms from $\llbracket \mathbf{Pth}^{(0,2)}_{\boldsymbol{\mathcal{A}}^{(2)}}
\rrbracket$ to $\mathbf{T}_{\Sigma}(X)$.
\end{restatable}

\begin{restatable}{proposition}{PDVDZIp}
\label{PDVDZIp} The many-sorted mapping $\mathrm{ip}^{(\llbracket 2\rrbracket,0)\sharp}$  is a $\Sigma$-homomorphism from  $\mathbf{T}_{\Sigma}(X)$ to $\llbracket \mathbf{Pth}^{(0,2)}_{\boldsymbol{\mathcal{A}}^{(2)}}
\rrbracket$.
\end{restatable}

\section{
\texorpdfstring
{A structure of $S$-sorted $2$-categorial $\Sigma$-algebra on $\llbracket \mathrm{Pth}_{\boldsymbol{\mathcal{A}}^{(2)}}\rrbracket$}
{A 2-categorial algebra on the second-order quotient}
} 

In this section we study the algebraic and categorial structure of the path quotient by the supremum of kernel of the second-order Curry-Howard mapping and the  congruence $\Upsilon^{[1]}$.

We now present the equations that satisfy the partial many-sorted $\Sigma^{\boldsymbol{\mathcal{A}}^{(2)}}$-algebra $ \llbracket\mathbf{Pth}_{\boldsymbol{\mathcal{A}}^{(2)}} \rrbracket $.

\begin{proposition}\label{PDVVarA2} Let $s$ be a sort in $S$ and $\llbracket \mathfrak{P}^{(2)}\rrbracket_{s}$ a second-order path class in $\llbracket\mathrm{Pth}_{\boldsymbol{\mathcal{A}}^{(2)}}\rrbracket_{s}$ then the following equalities hold
\allowdisplaybreaks
\begin{align*}
\mathrm{sc}^{0\llbracket\mathbf{Pth}_{\boldsymbol{\mathcal{A}}^{(2)}}\rrbracket}_{s}\left(
\mathrm{sc}^{0\llbracket\mathbf{Pth}_{\boldsymbol{\mathcal{A}}^{(2)}}\rrbracket}_{s}\left(
\llbracket
\mathfrak{P}^{(2)}
\rrbracket_{s}
\right)
\right)
& 
=
\mathrm{sc}^{0\llbracket\mathbf{Pth}_{\boldsymbol{\mathcal{A}}^{(2)}}\rrbracket}_{s}\left(
\llbracket
\mathfrak{P}^{(2)}
\rrbracket_{s}
\right);
\\
\mathrm{sc}^{0\llbracket\mathbf{Pth}_{\boldsymbol{\mathcal{A}}^{(2)}}\rrbracket}_{s}\left(
\mathrm{tg}^{0\llbracket\mathbf{Pth}_{\boldsymbol{\mathcal{A}}^{(2)}}\rrbracket}_{s}\left(
\llbracket
\mathfrak{P}^{(2)}
\rrbracket_{s}
\right)
\right)
& 
=
\mathrm{tg}^{0\llbracket\mathbf{Pth}_{\boldsymbol{\mathcal{A}}^{(2)}}\rrbracket}_{s}\left(
\llbracket
\mathfrak{P}^{(2)}
\rrbracket_{s}
\right);
\\
\mathrm{tg}^{0\llbracket\mathbf{Pth}_{\boldsymbol{\mathcal{A}}^{(2)}}\rrbracket}_{s}\left(
\mathrm{sc}^{0\llbracket\mathbf{Pth}_{\boldsymbol{\mathcal{A}}^{(2)}}\rrbracket}_{s}\left(
\llbracket
\mathfrak{P}^{(2)}
\rrbracket_{s}
\right)
\right)
& 
=
\mathrm{sc}^{0\llbracket\mathbf{Pth}_{\boldsymbol{\mathcal{A}}^{(2)}}\rrbracket}_{s}\left(
\llbracket
\mathfrak{P}^{(2)}
\rrbracket_{s}
\right);
\\
\mathrm{tg}^{0\llbracket\mathbf{Pth}_{\boldsymbol{\mathcal{A}}^{(2)}}\rrbracket}_{s}\left(
\mathrm{tg}^{0\llbracket\mathbf{Pth}_{\boldsymbol{\mathcal{A}}^{(2)}}\rrbracket}_{s}\left(
\llbracket
\mathfrak{P}^{(2)}
\rrbracket_{s}
\right)
\right)
& 
=
\mathrm{tg}^{0\llbracket\mathbf{Pth}_{\boldsymbol{\mathcal{A}}^{(2)}}\rrbracket}_{s}\left(
\llbracket
\mathfrak{P}^{(2)}
\rrbracket_{s}
\right).
\end{align*}	
\end{proposition}
\begin{proof}
We will only present the proof for the first equality. The other three are handled in a similar manner.

Note that the following chain of equalities holds
\allowdisplaybreaks
\begin{align*}
\mathrm{sc}^{0\mathbf{Pth}_{\boldsymbol{\mathcal{A}}^{(2)}}}_{s}\left(
\mathrm{sc}^{0\mathbf{Pth}_{\boldsymbol{\mathcal{A}}^{(2)}}}_{s}\left(
\mathfrak{P}^{(2)}
\right)\right)&=\mathrm{ip}^{(2,0)\sharp}_{s}\left(
\mathrm{sc}^{(0,2)}_{s}\left(
\mathrm{ip}^{(2,0)\sharp}_{s}\left(
\mathrm{sc}^{(0,2)}_{s}\left(
\mathfrak{P}^{(2)}
\right)\right)\right)\right)
\tag{1}
\\&=
\mathrm{ip}^{(2,0)\sharp}_{s}\left(
\mathrm{sc}^{(0,2)}_{s}\left(
\mathfrak{P}^{(2)}
\right)\right)
\tag{2}
\\&=
\mathrm{sc}^{0\mathbf{Pth}_{\boldsymbol{\mathcal{A}}^{(2)}}}_{s}\left(
\mathfrak{P}^{(2)}
\right).
\tag{3}
\end{align*}

In the just stated chain of equalities, the first equality follows from the definition of the $0$-source operation symbol in the many-sorted partial $\Sigma^{\boldsymbol{\mathcal{A}}^{(2)}}$-algebra $\mathbf{Pth}_{\boldsymbol{\mathcal{A}}^{(2)}}$ according to Proposition~\ref{PDPthCatAlg}; the second equality follows from Proposition~\ref{PDBasicEqZ}; finally, the last equality recovers from the definition of the $0$-source operation symbol in the many-sorted partial $\Sigma^{\boldsymbol{\mathcal{A}}^{(2)}}$-algebra $\mathbf{Pth}_{\boldsymbol{\mathcal{A}}^{(2)}}$ according to Proposition~\ref{PDPthCatAlg}.

Note that the desired equality follows from the fact that the $\llbracket\cdot\rrbracket$-classes of the two second-order paths above are equal.

This completes the proof.
\end{proof}

\begin{proposition}\label{PDVVarA3} Let $s$ be a sort in $S$ and $\llbracket \mathfrak{P}^{(2)}\rrbracket_{s}$, $\llbracket \mathfrak{Q}^{(2)}\rrbracket_{s}$  second-order path classes in $\llbracket\mathrm{Pth}_{\boldsymbol{\mathcal{A}}^{(2)}}\rrbracket_{s}$ then the following statements are equivalent 
\begin{enumerate}
\item[(i)] $\llbracket \mathfrak{Q}^{(2)}\rrbracket_{s}\circ^{0\llbracket\mathbf{Pth}_{\boldsymbol{\mathcal{A}}^{(2)}}\rrbracket}_{s} \llbracket \mathfrak{P}^{(2)}\rrbracket_{s}$ is defined;
\item[(ii)] $\mathrm{sc}^{0\llbracket\mathbf{Pth}_{\boldsymbol{\mathcal{A}}^{(2)}}\rrbracket}_{s}(
\llbracket\mathfrak{Q}^{(2)}\rrbracket_{s} ) = \mathrm{tg}^{0\llbracket\mathbf{Pth}_{\boldsymbol{\mathcal{A}}^{(2)}}\rrbracket}_{s}(
\llbracket\mathfrak{P}^{(2)}\rrbracket_{s} )$.
\end{enumerate}
\end{proposition}
\begin{proof}
To simplify the calculations we consider that, for some path terms $P,Q\in\mathrm{PT}_{\boldsymbol{\mathcal{A}},s}$, we have that 
\allowdisplaybreaks
\begin{align*}
\mathrm{sc}^{([1],2)}_{s}\left(\mathfrak{Q}^{(2)}\right)&=[Q]_{s};
&
\mathrm{tg}^{([1],2)}_{s}\left(\mathfrak{P}^{(2)}\right)&=[P]_{s}.
\end{align*}

The following chain of equivalences holds
\begin{flushleft}
$\llbracket \mathfrak{Q}^{(2)}\rrbracket_{s}\circ^{0\llbracket\mathbf{Pth}_{\boldsymbol{\mathcal{A}}^{(2)}}\rrbracket} \llbracket \mathfrak{P}^{(2)}\rrbracket_{s}$ is defined
\allowdisplaybreaks
\begin{align*}
\Leftrightarrow\quad& 
\mathfrak{Q}^{(2)} 
\circ^{0\mathbf{Pth}_{\boldsymbol{\mathcal{A}}^{(2)}}}_{s} 
\mathfrak{P}^{(2)}
\mbox{ is defined}
\tag{1}
\\
\Leftrightarrow\quad& 
\mathrm{sc}^{(0,2)}_{s}\left(
\mathfrak{Q}^{(2)}
\right)
=
\mathrm{tg}^{(0,2)}_{s}\left(
\mathfrak{P}^{(2)}
\right)
\tag{2}
\\
\Leftrightarrow\quad& 
\mathrm{sc}^{(0,[1])}_{s}\left(
\mathrm{ip}^{([1],X)@}_{s}\left(
[Q]_{s}
\right)
\right)
=
\mathrm{tg}^{(0,[1])}_{s}\left(
\mathrm{ip}^{([1],X)@}_{s}\left(
[P]_{s}
\right)
\right)
\tag{3}
\\
\Leftrightarrow\quad& 
\mathrm{sc}^{(0,[1])}_{s}\left(
\left[
\mathrm{ip}^{(1,X)@}_{s}(Q)
\right]_{s}
\right)
=
\mathrm{tg}^{(0,[1])}_{s}\left(
\left[
\mathrm{ip}^{(1,X)@}_{s}(P)
\right]_{s}
\right)
\tag{4}
\\
\Leftrightarrow\quad& 
\mathrm{sc}^{(0,1)}_{s}\left(
\mathrm{ip}^{(1,X)@}_{s}(Q)
\right)
=
\mathrm{tg}^{(0,1)}_{s}\left(
\mathrm{ip}^{(1,X)@}_{s}(P)
\right)
\tag{5}
\\
\Leftrightarrow\quad& 
\mathrm{sc}^{0[\mathbf{PT}_{\boldsymbol{\mathcal{A}}}]}_{s}\left(
[Q]_{s}
\right)
=
\mathrm{tg}^{0[\mathbf{PT}_{\boldsymbol{\mathcal{A}}}]}_{s}\left(
[P]_{s}
\right)
\tag{6}
\\
\Leftrightarrow\quad& 
\mathrm{ip}^{(2,[1])\sharp}_{s}\left(
\mathrm{sc}^{0[\mathbf{PT}_{\boldsymbol{\mathcal{A}}}]}_{s}\left(
[Q]_{s}
\right)\right)
=
\mathrm{ip}^{(2,[1])\sharp}_{s}\left(
\mathrm{tg}^{0[\mathbf{PT}_{\boldsymbol{\mathcal{A}}}]}_{s}\left(
[P]_{s}
\right)\right)
\tag{7}
\\
\Leftrightarrow\quad& 
\biggl\llbracket
\mathrm{ip}^{(2,[1])\sharp}_{s}\left(
\mathrm{sc}^{0[\mathbf{PT}_{\boldsymbol{\mathcal{A}}}]}_{s}\left(
[Q]_{s}
\right)\right)
\biggr\rrbracket_{s}
=
\biggl\llbracket
\mathrm{ip}^{(2,[1])\sharp}_{s}\left(
\mathrm{tg}^{0[\mathbf{PT}_{\boldsymbol{\mathcal{A}}}]}_{s}\left(
[P]_{s}
\right)\right)
\biggr\rrbracket_{s}
\tag{8}
\\
\Leftrightarrow\quad& 
\biggl\llbracket
\mathrm{sc}^{0\mathbf{Pth}_{\boldsymbol{\mathcal{A}}^{(2)}}}_{s}\left(
\mathfrak{Q}^{(2)}
\right)
\biggr\rrbracket_{s}
= 
\biggl\llbracket
\mathrm{tg}^{0\mathbf{Pth}_{\boldsymbol{\mathcal{A}}^{(2)}}}_{s}\left(
\mathfrak{P}^{(2)}
\right)
\biggr\rrbracket_{s}
\tag{9}
\\
\Leftrightarrow\quad& 
\mathrm{sc}^{0\llbracket\mathbf{Pth}_{\boldsymbol{\mathcal{A}}^{(2)}}\rrbracket}_{s}(
\llbracket\mathfrak{Q}^{(2)}\rrbracket_{s} ) = \mathrm{tg}^{0\llbracket\mathbf{Pth}_{\boldsymbol{\mathcal{A}}^{(2)}}\rrbracket}_{s}(
\llbracket\mathfrak{P}^{(2)}\rrbracket_{s} ).
\tag{10}
\end{align*}
\end{flushleft}

In the just stated chain of equivalences, the first equivalence follows from the interpretation of the $0$-composition operation in $\llbracket \mathbf{Pth}_{\boldsymbol{\mathcal{A}}^{(2)}}\rrbracket$ according to Proposition~\ref{PDVDCatAlg}; the second equivalence follows from Proposition~\ref{PDPthCatAlg}; the third equivalence follows from Definition~\ref{DDScTgZ}. Let us note  that we are using the abbreviations we introduced at the beginning of the proof; the fourth equivalence follows from Definition~\ref{DPTQIp}; the fifth equivalence follows from Proposition~\ref{PCHBasicEq}; the sixth equivalence follows from Proposition~\ref{PPTQVarA3}; the seventh equivalence follows from Proposition~\ref{PDBasicEq}; the eighth equivalence follows from Proposition~\ref{PDVDUId}; the ninth equivalence follows from Claim~\ref{CDPthCatAlgScZ}; finally, the last equivalence follows from the description of the $0$-source and $1$-source operations in $\llbracket \mathbf{Pth}_{\boldsymbol{\mathcal{A}}^{(2)}}\rrbracket$ according to Proposition~\ref{PDVDCatAlg}.
\end{proof}

\begin{proposition}\label{PDVVarA4} Let $s$ be a sort in $S$ and $\llbracket \mathfrak{P}^{(2)} \rrbracket_{s}$, $\llbracket \mathfrak{Q}^{(2)} \rrbracket_{s}$ second-order path classes in $\llbracket \mathrm{Pth}_{\boldsymbol{\mathcal{A}}^{(2)}}\rrbracket_{s}$. If its $0$-composition is defined, i.e., if 
\[
\mathrm{sc}^{0\llbracket\mathbf{Pth}_{\boldsymbol{\mathcal{A}}^{(2)}}\rrbracket}_{s}\left(
\mathfrak{Q}^{(2)}
\right)=
\mathrm{tg}^{0\llbracket\mathbf{Pth}_{\boldsymbol{\mathcal{A}}^{(2)}}\rrbracket}_{s}\left(
\mathfrak{P}^{(2)}\right),
\]
then the following equalities hold
\allowdisplaybreaks
\begin{align*}
\mathrm{sc}^{0\llbracket\mathbf{Pth}_{\boldsymbol{\mathcal{A}}^{(2)}}\rrbracket}_{s}\left(
\llbracket\mathfrak{Q}^{(2)}\rrbracket_{s} 
\circ^{0\llbracket\mathbf{Pth}_{\boldsymbol{\mathcal{A}}^{(2)}}\rrbracket}_{s}
\llbracket\mathfrak{P}^{(2)}\rrbracket_{s} 
\right)
&=
\mathrm{sc}^{0\llbracket\mathbf{Pth}_{\boldsymbol{\mathcal{A}}^{(2)}}\rrbracket}_{s}\left(
\llbracket\mathfrak{P}^{(2)}\rrbracket_{s} 
\right);
\\
\mathrm{tg}^{0\llbracket\mathbf{Pth}_{\boldsymbol{\mathcal{A}}^{(2)}}\rrbracket}_{s}\left(
\llbracket\mathfrak{Q}^{(2)}\rrbracket_{s} 
\circ^{0\llbracket\mathbf{Pth}_{\boldsymbol{\mathcal{A}}^{(2)}}\rrbracket}_{s}
\llbracket\mathfrak{P}^{(2)}\rrbracket_{s} 
\right)
&=
\mathrm{tg}^{0\llbracket\mathbf{Pth}_{\boldsymbol{\mathcal{A}}^{(2)}}\rrbracket}_{s}\left(
\llbracket\mathfrak{Q}^{(2)}\rrbracket_{s} 
\right).
\end{align*}
\end{proposition}
\begin{proof}
We will only present the proof for the first equality. The other one is handled in a similar manner.

The following chain of equalities holds
\begin{flushleft}
$\mathrm{sc}^{0\mathbf{Pth}_{\boldsymbol{\mathcal{A}}^{(2)}}}_{s}\left(
\mathfrak{Q}^{(2)}
\circ^{0\mathbf{Pth}_{\boldsymbol{\mathcal{A}}^{(2)}}}_{s}
\mathfrak{P}^{(2)}
\right)$
\allowdisplaybreaks
\begin{align*}
&=
\mathrm{ip}^{(2,[1])\sharp}_{s}\left(
\mathrm{sc}^{0[\mathbf{PT}_{\boldsymbol{\mathcal{A}}}]}_{s}\left(
\mathrm{sc}^{([1],2)}_{s}\left(
\mathfrak{Q}^{(2)}
\circ^{0\mathbf{Pth}_{\boldsymbol{\mathcal{A}}^{(2)}}}_{s}
\mathfrak{P}^{(2)}
\right)
\right)
\right)
\tag{1}
\\&=
\mathrm{ip}^{(2,[1])\sharp}_{s}\left(
\mathrm{sc}^{0[\mathbf{PT}_{\boldsymbol{\mathcal{A}}}]}_{s}\left(
\mathrm{sc}^{([1],2)}_{s}\left(
\mathfrak{Q}^{(2)}
\right)
\circ^{0[\mathbf{PT}_{\boldsymbol{\mathcal{A}}}]}_{s}
\mathrm{sc}^{([1],2)}_{s}\left(
\mathfrak{P}^{(2)}
\right)
\right)
\right)
\tag{2}
\\&=
\mathrm{ip}^{(2,[1])\sharp}_{s}\left(
\mathrm{sc}^{0[\mathbf{PT}_{\boldsymbol{\mathcal{A}}}]}_{s}\left(
\mathrm{sc}^{([1],2)}_{s}\left(
\mathfrak{P}^{(2)}
\right)
\right)
\right)
\tag{3}
\\&=
\mathrm{sc}^{0\mathbf{Pth}_{\boldsymbol{\mathcal{A}}^{(2)}}}_{s}\left(
\mathfrak{P}^{(2)}
\right).
\tag{4}
\end{align*}
\end{flushleft}

In the just stated chain of equalities, the first equality follows from Claim~\ref{CDPthCatAlgScZ}; the second equality follows from Claim~\ref{CDPthCatAlgCompZ}; the third equality follows from Proposition~\ref{PPTQVarA4}; finally, the last equality follows from Claim~\ref{CDPthCatAlgScZ}.

Note that the desired equality follows from the fact that the $\llbracket\cdot\rrbracket$-classes of the two second-order paths above are equal.

This completes the proof.
\end{proof}

\begin{proposition}\label{PDVVarA5} Let $s$ be a sort in $S$ and $\llbracket \mathfrak{P}^{(2)} \rrbracket_{s}$ a second-order path class in $\llbracket \mathrm{Pth}_{\boldsymbol{\mathcal{A}}^{(2)}}\rrbracket_{s}$, then the following equalities hold
\allowdisplaybreaks
\begin{align*}
\llbracket\mathfrak{P}^{(2)}\rrbracket_{s} 
\circ^{0\llbracket\mathbf{Pth}_{\boldsymbol{\mathcal{A}}^{(2)}}\rrbracket}_{s}
\mathrm{sc}^{0\llbracket\mathbf{Pth}_{\boldsymbol{\mathcal{A}}^{(2)}}\rrbracket}_{s}\left(
\llbracket\mathfrak{P}^{(2)}\rrbracket_{s} 
\right)
&=
\llbracket\mathfrak{P}^{(2)}\rrbracket_{s};
\\
\mathrm{tg}^{0\llbracket\mathbf{Pth}_{\boldsymbol{\mathcal{A}}^{(2)}}\rrbracket}_{s}\left(
\llbracket\mathfrak{P}^{(2)}\rrbracket_{s} 
\right)
\circ^{0\llbracket\mathbf{Pth}_{\boldsymbol{\mathcal{A}}^{(2)}}\rrbracket}_{s}
\llbracket\mathfrak{P}^{(2)}\rrbracket_{s} 
&=
\llbracket\mathfrak{P}^{(2)}\rrbracket_{s};
\end{align*}
\end{proposition}
\begin{proof}
We will only present the proof for the first equality. The other one is handled in a similar manner.

Let us recall from Definitions~\ref{DDUps} and~\ref{DDUpsCong} that the following equality holds 
\[
\left[
\mathfrak{P}^{(2)}\circ^{0\mathbf{Pth}_{\boldsymbol{\mathcal{A}}^{(2)}}}
\mathrm{sc}^{0\mathbf{Pth}_{\boldsymbol{\mathcal{A}}^{(2)}}}_{s}\left(
\mathfrak{P}^{(2)}
\right)
\right]_{\Upsilon^{[1]}_{s}}
=
\left[
\mathfrak{P}^{(2)}
\right]_{\Upsilon^{[1]}_{s}}
.
\]

Therefore, taking $\llbracket \cdot \rrbracket$-classes, we conclude that 
\[
\Bigl\llbracket
\mathfrak{P}^{(2)}\circ^{0\mathbf{Pth}_{\boldsymbol{\mathcal{A}}^{(2)}}}
\mathrm{sc}^{0\mathbf{Pth}_{\boldsymbol{\mathcal{A}}^{(2)}}}_{s}\left(
\mathfrak{P}^{(2)}
\right)
\Bigr\rrbracket_{s}
=
\Bigl\llbracket
\mathfrak{P}^{(2)}
\Bigr\rrbracket_{s}
.
\]

Taking into account the definition of the operation symbols in $\llbracket \mathbf{Pth}_{\boldsymbol{\mathcal{A}}^{(2)}} \rrbracket$, according to Proposition~\ref{PDVDCatAlg}, we conclude that 
\[
\Bigl\llbracket\mathfrak{P}^{(2)}\Bigr\rrbracket_{s} 
\circ^{0\llbracket\mathbf{Pth}_{\boldsymbol{\mathcal{A}}^{(2)}}\rrbracket}_{s}
\mathrm{sc}^{0\llbracket\mathbf{Pth}_{\boldsymbol{\mathcal{A}}^{(2)}}\rrbracket}_{s}\left(
\Bigl\llbracket\mathfrak{P}^{(2)}\Bigr\rrbracket_{s} 
\right)
=
\Bigl\llbracket\mathfrak{P}^{(2)}\Bigr\rrbracket_{s}.
\]

This completes the proof.
\end{proof}

\begin{proposition}\label{PDVVarA6} Let $s$ be a sort in $S$ and $\llbracket \mathfrak{P}^{(2)} \rrbracket_{s}$, $\llbracket \mathfrak{Q}^{(2)} \rrbracket_{s}$, $\llbracket \mathfrak{R}^{(2)} \rrbracket_{s}$  second-order path classes in $\llbracket \mathrm{Pth}_{\boldsymbol{\mathcal{A}}^{(2)}}\rrbracket_{s}$ satisfying that 
\allowdisplaybreaks
\begin{align*}
\mathrm{sc}^{0\llbracket \mathbf{Pth}_{\boldsymbol{\mathcal{A}}^{(2)}} \rrbracket}_{s}\left(
\mathfrak{R}^{(2)}
\right)
&=
\mathrm{tg}^{0\llbracket \mathbf{Pth}_{\boldsymbol{\mathcal{A}}^{(2)}} \rrbracket}_{s}\left(
\mathfrak{Q}^{(2)}
\right);
\\
\mathrm{sc}^{0\llbracket \mathbf{Pth}_{\boldsymbol{\mathcal{A}}^{(2)}} \rrbracket}_{s}\left(
\mathfrak{Q}^{(2)}
\right)
&=
\mathrm{tg}^{0\llbracket \mathbf{Pth}_{\boldsymbol{\mathcal{A}}^{(2)}} \rrbracket}_{s}\left(
\mathfrak{P}^{(2)}
\right),
\end{align*}
then the following equality holds
\allowdisplaybreaks
\begin{multline*}
\llbracket
\mathfrak{R}^{(2)}
\rrbracket_{s}
\circ^{0\llbracket \mathbf{Pth}_{\boldsymbol{\mathcal{A}}^{(2)}}\rrbracket}_{s}
\left(
\llbracket
\mathfrak{Q}^{(2)}
\rrbracket_{s}
\circ^{0\llbracket \mathbf{Pth}_{\boldsymbol{\mathcal{A}}^{(2)}}\rrbracket}_{s}
\llbracket
\mathfrak{P}^{(2)}
\rrbracket_{s}
\right)
\\=
\left(
\llbracket
\mathfrak{R}^{(2)}
\rrbracket_{s}
\circ^{0\llbracket \mathbf{Pth}_{\boldsymbol{\mathcal{A}}^{(2)}}\rrbracket}_{s}
\llbracket
\mathfrak{Q}^{(2)}
\rrbracket_{s}
\right)
\circ^{0\llbracket \mathbf{Pth}_{\boldsymbol{\mathcal{A}}^{(2)}}\rrbracket}_{s}
\llbracket
\mathfrak{P}^{(2)}
\rrbracket_{s}.
\end{multline*}
\end{proposition}
\begin{proof}
Let us recall from Definitions~\ref{DDUps} and~\ref{DDUpsCong} that the following equality holds 
\[
\left[
\mathfrak{R}^{2}
\circ^{0\mathbf{Pth}_{\boldsymbol{\mathcal{A}}^{(2)}}}_{s}
\left(
\mathfrak{Q}^{2}
\circ^{0\mathbf{Pth}_{\boldsymbol{\mathcal{A}}^{(2)}}}_{s}
\mathfrak{P}^{2}
\right)
\right]_{\Upsilon^{[1]}_{s}}
=
\left[
\left(
\mathfrak{R}^{2}
\circ^{0\mathbf{Pth}_{\boldsymbol{\mathcal{A}}^{(2)}}}_{s}
\mathfrak{Q}^{2}
\right)
\circ^{0\mathbf{Pth}_{\boldsymbol{\mathcal{A}}^{(2)}}}_{s}
\mathfrak{P}^{2}
\right]_{\Upsilon^{[1]}_{s}}
.
\]

Therefore, taking $\llbracket \cdot \rrbracket$-classes, we conclude that 
\[
\Bigl\llbracket
\mathfrak{R}^{2}
\circ^{0\mathbf{Pth}_{\boldsymbol{\mathcal{A}}^{(2)}}}_{s}
\left(
\mathfrak{Q}^{2}
\circ^{0\mathbf{Pth}_{\boldsymbol{\mathcal{A}}^{(2)}}}_{s}
\mathfrak{P}^{2}
\right)
\Bigr\rrbracket_{s}
=
\Bigl\llbracket
\left(
\mathfrak{R}^{2}
\circ^{0\mathbf{Pth}_{\boldsymbol{\mathcal{A}}^{(2)}}}_{s}
\mathfrak{Q}^{2}
\right)
\circ^{0\mathbf{Pth}_{\boldsymbol{\mathcal{A}}^{(2)}}}_{s}
\mathfrak{P}^{2}
\Bigr\rrbracket_{s}
.
\]

Taking into account the definition of the operation symbols in $\llbracket \mathbf{Pth}_{\boldsymbol{\mathcal{A}}^{(2)}} \rrbracket$, according to Proposition~\ref{PDVDCatAlg}, we conclude that 
\allowdisplaybreaks
\begin{multline*}
\llbracket
\mathfrak{R}^{(2)}
\rrbracket_{s}
\circ^{0\llbracket \mathbf{Pth}_{\boldsymbol{\mathcal{A}}^{(2)}}\rrbracket}_{s}
\left(
\llbracket
\mathfrak{Q}^{(2)}
\rrbracket_{s}
\circ^{0\llbracket \mathbf{Pth}_{\boldsymbol{\mathcal{A}}^{(2)}}\rrbracket}_{s}
\llbracket
\mathfrak{P}^{(2)}
\rrbracket_{s}
\right)
\\=
\left(
\llbracket
\mathfrak{R}^{(2)}
\rrbracket_{s}
\circ^{0\llbracket \mathbf{Pth}_{\boldsymbol{\mathcal{A}}^{(2)}}\rrbracket}_{s}
\llbracket
\mathfrak{Q}^{(2)}
\rrbracket_{s}
\right)
\circ^{0\llbracket \mathbf{Pth}_{\boldsymbol{\mathcal{A}}^{(2)}}\rrbracket}_{s}
\llbracket
\mathfrak{P}^{(2)}
\rrbracket_{s}.
\end{multline*}

This completes the proof.
\end{proof}

\begin{proposition}\label{PDVVarB2} Let $s$ be a sort in $S$ and $\llbracket \mathfrak{P}^{(2)}\rrbracket_{s}$ a second-order path class in $\llbracket\mathrm{Pth}_{\boldsymbol{\mathcal{A}}^{(2)}}\rrbracket_{s}$ then the following equalities holds
\allowdisplaybreaks
\begin{align*}
\mathrm{sc}^{1\llbracket\mathbf{Pth}_{\boldsymbol{\mathcal{A}}^{(2)}}\rrbracket}_{s}\left(
\mathrm{sc}^{1\llbracket\mathbf{Pth}_{\boldsymbol{\mathcal{A}}^{(2)}}\rrbracket}_{s}\left(
\llbracket
\mathfrak{P}^{(2)}
\rrbracket_{s}
\right)
\right)
&
=
\mathrm{sc}^{1\llbracket\mathbf{Pth}_{\boldsymbol{\mathcal{A}}^{(2)}}\rrbracket}_{s}\left(
\llbracket
\mathfrak{P}^{(2)}
\rrbracket_{s}
\right);
\\
\mathrm{sc}^{1\llbracket\mathbf{Pth}_{\boldsymbol{\mathcal{A}}^{(2)}}\rrbracket}_{s}\left(
\mathrm{tg}^{1\llbracket\mathbf{Pth}_{\boldsymbol{\mathcal{A}}^{(2)}}\rrbracket}_{s}\left(
\llbracket
\mathfrak{P}^{(2)}
\rrbracket_{s}
\right)
\right)
&
=
\mathrm{tg}^{1\llbracket\mathbf{Pth}_{\boldsymbol{\mathcal{A}}^{(2)}}\rrbracket}_{s}\left(
\llbracket
\mathfrak{P}^{(2)}
\rrbracket_{s}
\right);
\\
\mathrm{tg}^{1\llbracket\mathbf{Pth}_{\boldsymbol{\mathcal{A}}^{(2)}}\rrbracket}_{s}\left(
\mathrm{sc}^{1\llbracket\mathbf{Pth}_{\boldsymbol{\mathcal{A}}^{(2)}}\rrbracket}_{s}\left(
\llbracket
\mathfrak{P}^{(2)}
\rrbracket_{s}
\right)
\right)
&
=
\mathrm{sc}^{1\llbracket\mathbf{Pth}_{\boldsymbol{\mathcal{A}}^{(2)}}\rrbracket}_{s}\left(
\llbracket
\mathfrak{P}^{(2)}
\rrbracket_{s}
\right);
\\
\mathrm{tg}^{1\llbracket\mathbf{Pth}_{\boldsymbol{\mathcal{A}}^{(2)}}\rrbracket}_{s}\left(
\mathrm{tg}^{1\llbracket\mathbf{Pth}_{\boldsymbol{\mathcal{A}}^{(2)}}\rrbracket}_{s}\left(
\llbracket
\mathfrak{P}^{(2)}
\rrbracket_{s}
\right)
\right)
&
=
\mathrm{tg}^{1\llbracket\mathbf{Pth}_{\boldsymbol{\mathcal{A}}^{(2)}}\rrbracket}_{s}\left(
\llbracket
\mathfrak{P}^{(2)}
\rrbracket_{s}
\right).
\end{align*}
\end{proposition}
\begin{proof}
We will only present the proof for the first equality. The other three are handled in a similar manner.

The following chain of equalities holds
\begin{flushleft}
$\mathrm{sc}^{1\mathbf{Pth}_{\boldsymbol{\mathcal{A}}^{(2)}}}_{s}\left(
\mathrm{sc}^{1\mathbf{Pth}_{\boldsymbol{\mathcal{A}}^{(2)}}}_{s}\left(
\mathfrak{P}^{(2)}
\right)\right)$
\allowdisplaybreaks
\begin{align*}
\qquad&=
\mathrm{ip}^{(2,[1])\sharp}_{s}\left(
\mathrm{sc}^{([1],2)}_{s}\left(
\mathrm{ip}^{(2,[1])\sharp}_{s}\left(
\mathrm{sc}^{([1],2)}_{s}\left(
\mathfrak{P}^{(2)}
\right)\right)\right)\right)
\tag{1}
\\&=
\mathrm{ip}^{(2,[1])\sharp}_{s}\left(
\mathrm{sc}^{([1],2)}_{s}\left(
\mathfrak{P}^{(2)}
\right)\right)
\tag{2}
\\&=
\mathrm{sc}^{1\mathbf{Pth}_{\boldsymbol{\mathcal{A}}^{(2)}}}_{s}\left(
\mathfrak{P}^{(2)}
\right).
\tag{3}
\end{align*}
\end{flushleft}

In the just stated chain of equalities, the first equality unravels the description of the $1$-source operation symbol in the many-sorted partial $\Sigma^{\boldsymbol{\mathcal{A}}^{(2)}}$-algebra $\mathbf{Pth}_{\boldsymbol{\mathcal{A}}^{(2)}}$ according to Proposition~\ref{PDPthDCatAlg}; the second equality follows from Proposition~\ref{PDBasicEq}; finally, the last equivalence recovers the description of the $1$-source operation symbol in the many-sorted partial $\Sigma^{\boldsymbol{\mathcal{A}}^{(2)}}$-algebra $\mathbf{Pth}_{\boldsymbol{\mathcal{A}}^{(2)}}$ according to Proposition~\ref{PDPthDCatAlg}.

Note that the desired equality follows from the fact that the $\llbracket\cdot\rrbracket$-classes of the two second-order paths above are equal.

This completes the proof.
\end{proof}

\begin{proposition}\label{PDVVarB3} Let $s$ be a sort in $S$ and $\llbracket \mathfrak{P}^{(2)}\rrbracket_{s}$, $\llbracket \mathfrak{Q}^{(2)}\rrbracket_{s}$  second-order path classes in $\llbracket\mathrm{Pth}_{\boldsymbol{\mathcal{A}}^{(2)}}\rrbracket_{s}$ then the following statements are equivalent 
\begin{enumerate}
\item[(i)] $\llbracket \mathfrak{Q}^{(2)}\rrbracket_{s}\circ^{1\llbracket\mathbf{Pth}_{\boldsymbol{\mathcal{A}}^{(2)}}\rrbracket} \llbracket \mathfrak{P}^{(2)}\rrbracket_{s}$ is defined;
\item[(ii)] $\mathrm{sc}^{1\llbracket\mathbf{Pth}_{\boldsymbol{\mathcal{A}}^{(2)}}\rrbracket}_{s}(
\llbracket\mathfrak{Q}^{(2)}\rrbracket_{s} ) = \mathrm{tg}^{1\llbracket\mathbf{Pth}_{\boldsymbol{\mathcal{A}}^{(2)}}\rrbracket}_{s}(
\llbracket\mathfrak{P}^{(2)}\rrbracket_{s} )$.
\end{enumerate}
\end{proposition}
\begin{proof}
The following chain of equivalences holds
\begin{flushleft}
$\llbracket \mathfrak{Q}^{(2)}\rrbracket_{s}\circ^{1\llbracket\mathbf{Pth}_{\boldsymbol{\mathcal{A}}^{(2)}}\rrbracket} \llbracket \mathfrak{P}^{(2)}\rrbracket_{s}$ is defined
\allowdisplaybreaks
\begin{align*}
\qquad&\Leftrightarrow\quad
 \mathfrak{Q}^{(2)}\circ^{1\mathbf{Pth}_{\boldsymbol{\mathcal{A}}^{(2)}}}  \mathfrak{P}^{(2)} \mbox{ is defined}
 \tag{1}
 \\
 \qquad&\Leftrightarrow\quad
 \mathrm{sc}^{([1],2)}_{s}\left(
 \mathfrak{Q}^{(2)}
 \right)
 =
 \mathrm{tg}^{([1],2)}_{s}\left(
 \mathfrak{P}^{(2)}
 \right) 
 \tag{2}
\\
 \qquad&\Leftrightarrow\quad
 \mathrm{ip}^{(2,[1])\sharp}_{s}\left(
 \mathrm{sc}^{([1],2)}_{s}\left(
 \mathfrak{Q}^{(2)}
 \right)\right)
 =
 \mathrm{ip}^{(2,[1])\sharp}_{s}\left(
 \mathrm{tg}^{([1],2)}_{s}\left(
 \mathfrak{P}^{(2)}
 \right)  \right) 
 \tag{3}
 \\
 \qquad&\Leftrightarrow\quad
 \Bigl\llbracket
 \mathrm{ip}^{(2,[1])\sharp}_{s}\left(
 \mathrm{sc}^{([1],2)}_{s}\left(
 \mathfrak{Q}^{(2)}
 \right)\right)
 \Bigr\rrbracket_{s}
 =
 \Bigl\llbracket
 \mathrm{ip}^{(2,[1])\sharp}_{s}\left(
 \mathrm{tg}^{([1],2)}_{s}\left(
 \mathfrak{P}^{(2)}
 \right)  \right) 
 \Bigr\rrbracket_{s}
 \tag{4}
 \\
 \qquad&\Leftrightarrow\quad
 \Bigl\llbracket
 \mathrm{sc}^{1\mathbf{Pth}_{\boldsymbol{\mathcal{A}}^{(2)}}}_{s}\left(
 \mathfrak{Q}^{(2)}
 \right)
 \Bigr\rrbracket_{s} 
 =
  \Bigl\llbracket
 \mathrm{tg}^{1\mathbf{Pth}_{\boldsymbol{\mathcal{A}}^{(2)}}}_{s}\left(
 \mathfrak{P}^{(2)}
 \right)
 \Bigr\rrbracket_{s} 
 \tag{5}
\\ 
 \qquad&\Leftrightarrow\quad
 \mathrm{sc}^{1\llbracket\mathbf{Pth}_{\boldsymbol{\mathcal{A}}^{(2)}}\rrbracket}_{s}
 \left(
\llbracket\mathfrak{Q}^{(2)}\rrbracket_{s} 
\right) = 
\mathrm{tg}^{1\llbracket\mathbf{Pth}_{\boldsymbol{\mathcal{A}}^{(2)}}\rrbracket}_{s}
\left(
\llbracket\mathfrak{P}^{(2)}\rrbracket_{s} 
\right).
\tag{6}
\end{align*}
\end{flushleft}

In the just stated chain of equivalences, the first equivalence follows from Proposition~\ref{PDVDCatAlg}; the second equivalence follows from Definition~\ref{DDPthComp}; the third equivalence follows from Proposition~\ref{PDBasicEq}; the fourth equivalence follows from Proposition~\ref{PDVDUId}; the fifth equivalence follows from the description of the $1$-source and $1$-target operation symbols in the many-sorted partial $\Sigma^{\boldsymbol{\mathcal{A}}^{(2)}}$-algebra $\mathbf{Pth}_{\boldsymbol{\mathcal{A}}^{(2)}}$, according to Proposition~\ref{PDPthDCatAlg}; finally, the last equality follows form description of the $1$-source and $1$-target operation symbols in the many-sorted partial $\Sigma^{\boldsymbol{\mathcal{A}}^{(2)}}$-algebra $\llbracket \mathbf{Pth}_{\boldsymbol{\mathcal{A}}^{(2)}}\rrbracket $, according to Proposition~\ref{PDVDCatAlg}.

This completes the proof.
\end{proof}

\begin{proposition}\label{PDVVarB4} Let $s$ be a sort in $S$ and $\llbracket \mathfrak{P}^{(2)} \rrbracket_{s}$, $\llbracket \mathfrak{Q}^{(2)} \rrbracket_{s}$ second-order path classes in $\llbracket \mathrm{Pth}_{\boldsymbol{\mathcal{A}}^{(2)}}\rrbracket_{s}$. If its $1$-composition is defined, i.e., if 
\[
\mathrm{sc}^{1\llbracket\mathbf{Pth}_{\boldsymbol{\mathcal{A}}^{(2)}}\rrbracket}_{s}\left(
\mathfrak{Q}^{(2)}\right)
=
\mathrm{tg}^{1\llbracket\mathbf{Pth}_{\boldsymbol{\mathcal{A}}^{(2)}}\rrbracket}_{s}\left(
\mathfrak{P}^{(2)}\right),
\]
the following equalities hold
\allowdisplaybreaks
\begin{align*}
\mathrm{sc}^{1\llbracket\mathbf{Pth}_{\boldsymbol{\mathcal{A}}^{(2)}}\rrbracket}_{s}\left(
\llbracket\mathfrak{Q}^{(2)}\rrbracket_{s} 
\circ^{1\llbracket\mathbf{Pth}_{\boldsymbol{\mathcal{A}}^{(2)}}\rrbracket}_{s}
\llbracket\mathfrak{P}^{(2)}\rrbracket_{s} 
\right)
&=
\mathrm{sc}^{1\llbracket\mathbf{Pth}_{\boldsymbol{\mathcal{A}}^{(2)}}\rrbracket}_{s}\left(
\llbracket\mathfrak{P}^{(2)}\rrbracket_{s} 
\right);
\\
\mathrm{tg}^{1\llbracket\mathbf{Pth}_{\boldsymbol{\mathcal{A}}^{(2)}}\rrbracket}_{s}\left(
\llbracket\mathfrak{Q}^{(2)}\rrbracket_{s} 
\circ^{1\llbracket\mathbf{Pth}_{\boldsymbol{\mathcal{A}}^{(2)}}\rrbracket}_{s}
\llbracket\mathfrak{P}^{(2)}\rrbracket_{s} 
\right)
&=
\mathrm{tg}^{1\llbracket\mathbf{Pth}_{\boldsymbol{\mathcal{A}}^{(2)}}\rrbracket}_{s}\left(
\llbracket\mathfrak{Q}^{(2)}\rrbracket_{s} 
\right).
\end{align*}
\end{proposition}
\begin{proof}
We will only present the proof for the first equality. The other one is handled in a similar manner.

The following chain of equalities holds.
\allowdisplaybreaks
\begin{align*}
\mathrm{sc}^{1\mathbf{Pth}_{\boldsymbol{\mathcal{A}}^{(2)}}}_{s}\left(
\mathfrak{Q}^{(2)} 
\circ^{1\mathbf{Pth}_{\boldsymbol{\mathcal{A}}^{(2)}}}_{s}
\mathfrak{P}^{(2)} 
\right)&=
\mathrm{ip}^{(2,[1])\sharp}_{s}\left(
\mathrm{sc}^{([1],2)}_{s}\left(
\mathfrak{Q}^{(2)} 
\circ^{1\mathbf{Pth}_{\boldsymbol{\mathcal{A}}^{(2)}}}_{s}
\mathfrak{P}^{(2)} 
\right)\right)
\tag{1}
\\&=
\mathrm{ip}^{(2,[1])\sharp}_{s}\left(
\mathrm{sc}^{([1],2)}_{s}\left(
\mathfrak{P}^{(2)} 
\right)\right)
\tag{2}
\\&=
\mathrm{sc}^{1\mathbf{Pth}_{\boldsymbol{\mathcal{A}}^{(2)}}}_{s}\left(
\mathfrak{P}^{(2)} 
\right).
\tag{3}
\end{align*}

In the just stated chain of equalities, the first equality unravels the description of the $1$-source  operation symbol in the many-sorted partial $\Sigma^{\boldsymbol{\mathcal{A}}^{(2)}}$-algebra $\mathbf{Pth}_{\boldsymbol{\mathcal{A}}^{(2)}}$, according to Proposition~\ref{PDPthDCatAlg}; the second equality follows from Definition~\ref{DDPthComp}; finally the last equality recovers the description of the $1$-source  operation symbol in the many-sorted partial $\Sigma^{\boldsymbol{\mathcal{A}}^{(2)}}$-algebra $\mathbf{Pth}_{\boldsymbol{\mathcal{A}}^{(2)}}$, according to Proposition~\ref{PDPthDCatAlg}.

Note that the desired equality follows from the fact that the $\llbracket\cdot\rrbracket$-classes of the two second-order paths above are equal.

This completes the proof.
\end{proof}

\begin{proposition}\label{PDVVarB5} Let $s$ be a sort in $S$ and $\llbracket \mathfrak{P}^{(2)} \rrbracket_{s}$ a second-order path class in $\llbracket \mathrm{Pth}_{\boldsymbol{\mathcal{A}}^{(2)}}\rrbracket_{s}$, then the following equalities hold
\allowdisplaybreaks
\begin{align*}
\llbracket\mathfrak{P}^{(2)}\rrbracket_{s} 
\circ^{1\llbracket\mathbf{Pth}_{\boldsymbol{\mathcal{A}}^{(2)}}\rrbracket}_{s}
\mathrm{sc}^{1\llbracket\mathbf{Pth}_{\boldsymbol{\mathcal{A}}^{(2)}}\rrbracket}_{s}\left(
\llbracket\mathfrak{P}^{(2)}\rrbracket_{s} 
\right)
&=
\llbracket\mathfrak{P}^{(2)}\rrbracket_{s};
\\
\mathrm{tg}^{1\llbracket\mathbf{Pth}_{\boldsymbol{\mathcal{A}}^{(2)}}\rrbracket}_{s}\left(
\llbracket\mathfrak{P}^{(2)}\rrbracket_{s} 
\right)
\circ^{1\llbracket\mathbf{Pth}_{\boldsymbol{\mathcal{A}}^{(2)}}\rrbracket}_{s}
\llbracket\mathfrak{P}^{(2)}\rrbracket_{s} 
&=
\llbracket\mathfrak{P}^{(2)}\rrbracket_{s}.
\end{align*}
\end{proposition}
\begin{proof}
We will only present the proof for the first equality. The other one is handled in a similar manner.

The following equality holds according to Proposition~\ref{PDPthComp}
\[
\mathfrak{P}^{(2)}
\circ^{1\mathbf{Pth}_{\boldsymbol{\mathcal{A}}^{(2)}}}_{s}
\mathrm{sc}^{1\mathbf{Pth}_{\boldsymbol{\mathcal{A}}^{(2)}}}_{s}\left(
\mathfrak{P}^{(2)}
\right)
=
\mathfrak{P}^{(2)}
.
\]

Note that the desired equality follows from the fact that the $\llbracket\cdot\rrbracket$-classes of the two second-order paths above are equal.

This completes the proof.
\end{proof}

\begin{proposition}\label{PDVVarB6} Let $s$ be a sort in $S$ and $\llbracket \mathfrak{P}^{(2)} \rrbracket_{s}$, $\llbracket \mathfrak{Q}^{(2)} \rrbracket_{s}$, $\llbracket \mathfrak{R}^{(2)} \rrbracket_{s}$  second-order path classes in $\llbracket \mathrm{Pth}_{\boldsymbol{\mathcal{A}}^{(2)}}\rrbracket_{s}$ satisfying that
\allowdisplaybreaks
\begin{align*}
\mathrm{sc}^{1\llbracket \mathbf{Pth}_{\boldsymbol{\mathcal{A}}^{(2)}}\rrbracket}_{s}\left(
\llbracket \mathfrak{R}^{(2)}\rrbracket_{s}
\right)
&=
\mathrm{tg}^{1\llbracket \mathbf{Pth}_{\boldsymbol{\mathcal{A}}^{(2)}}\rrbracket}_{s}\left(
\llbracket \mathfrak{Q}^{(2)}\rrbracket_{s}
\right)
\\
\mathrm{sc}^{1\llbracket \mathbf{Pth}_{\boldsymbol{\mathcal{A}}^{(2)}}\rrbracket}_{s}\left(
\llbracket \mathfrak{Q}^{(2)}\rrbracket_{s}
\right)
&=
\mathrm{tg}^{1\llbracket \mathbf{Pth}_{\boldsymbol{\mathcal{A}}^{(2)}}\rrbracket}_{s}\left(
\llbracket \mathfrak{P}^{(2)}\rrbracket_{s}
\right),
\end{align*}
then the following equalities hold
\allowdisplaybreaks
\begin{multline*}
\llbracket
\mathfrak{R}^{(2)}
\rrbracket_{s}
\circ^{1\llbracket \mathbf{Pth}_{\boldsymbol{\mathcal{A}}^{(2)}}\rrbracket}_{s}
\left(
\llbracket
\mathfrak{Q}^{(2)}
\rrbracket_{s}
\circ^{1\llbracket \mathbf{Pth}_{\boldsymbol{\mathcal{A}}^{(2)}}\rrbracket}_{s}
\llbracket
\mathfrak{P}^{(2)}
\rrbracket_{s}
\right)
\\=
\left(
\llbracket
\mathfrak{R}^{(2)}
\rrbracket_{s}
\circ^{1\llbracket \mathbf{Pth}_{\boldsymbol{\mathcal{A}}^{(2)}}\rrbracket}_{s}
\llbracket
\mathfrak{Q}^{(2)}
\rrbracket_{s}
\right)
\circ^{1\llbracket \mathbf{Pth}_{\boldsymbol{\mathcal{A}}^{(2)}}\rrbracket}_{s}
\llbracket
\mathfrak{P}^{(2)}
\rrbracket_{s}.
\end{multline*}
\end{proposition}
\begin{proof}
The following equality holds according to Proposition~\ref{PDPthComp}
\[
\mathfrak{R}^{(2)}
\circ^{1\mathbf{Pth}_{\boldsymbol{\mathcal{A}}^{(2)}}}_{s}
\left(
\mathfrak{Q}^{(2)}
\circ^{1\mathbf{Pth}_{\boldsymbol{\mathcal{A}}^{(2)}}}_{s}
\mathfrak{P}^{(2)}
\right)
=
\left(
\mathfrak{R}^{(2)}
\circ^{1\mathbf{Pth}_{\boldsymbol{\mathcal{A}}^{(2)}}}_{s}
\mathfrak{Q}^{(2)}
\right)
\circ^{1\mathbf{Pth}_{\boldsymbol{\mathcal{A}}^{(2)}}}_{s}
\mathfrak{P}^{(2)}
.
\]

Note that the desired equality follows from the fact that the $\llbracket\cdot\rrbracket$-classes of the two second-order paths above are equal.

This completes the proof.
\end{proof}

\begin{proposition}\label{PDVVarAB1} Let $s$ be a sort in $S$ and $\llbracket \mathfrak{P}^{(2)} \rrbracket_{s}$ a  second-order path class in $\llbracket \mathrm{Pth}_{\boldsymbol{\mathcal{A}}^{(2)}}\rrbracket_{s}$, then the following equalities hold
\allowdisplaybreaks
\begin{align*}
\mathrm{sc}^{1\llbracket \mathbf{Pth}_{\boldsymbol{\mathcal{A}}^{(2)}}\rrbracket}_{s}\left(
\mathrm{sc}^{0\llbracket \mathbf{Pth}_{\boldsymbol{\mathcal{A}}^{(2)}}\rrbracket}_{s}\left(
\llbracket \mathfrak{P}^{(2)}\rrbracket_{s}
\right)\right)
&=
\mathrm{sc}^{0\llbracket \mathbf{Pth}_{\boldsymbol{\mathcal{A}}^{(2)}}\rrbracket}_{s}\left(
\llbracket \mathfrak{P}^{(2)}\rrbracket_{s}
\right);
\\
\mathrm{sc}^{0\llbracket \mathbf{Pth}_{\boldsymbol{\mathcal{A}}^{(2)}}\rrbracket}_{s}\left(
\mathrm{sc}^{1\llbracket \mathbf{Pth}_{\boldsymbol{\mathcal{A}}^{(2)}}\rrbracket}_{s}\left(
\llbracket \mathfrak{P}^{(2)}\rrbracket_{s}
\right)\right)
&=
\mathrm{sc}^{0\llbracket \mathbf{Pth}_{\boldsymbol{\mathcal{A}}^{(2)}}\rrbracket}_{s}\left(
\llbracket \mathfrak{P}^{(2)}\rrbracket_{s}
\right);
\\
\mathrm{sc}^{0\llbracket \mathbf{Pth}_{\boldsymbol{\mathcal{A}}^{(2)}}\rrbracket}_{s}\left(
\mathrm{tg}^{1\llbracket \mathbf{Pth}_{\boldsymbol{\mathcal{A}}^{(2)}}\rrbracket}_{s}\left(
\llbracket \mathfrak{P}^{(2)}\rrbracket_{s}
\right)\right)
&=
\mathrm{sc}^{0\llbracket \mathbf{Pth}_{\boldsymbol{\mathcal{A}}^{(2)}}\rrbracket}_{s}\left(
\llbracket \mathfrak{P}^{(2)}\rrbracket_{s}
\right);
\\
\mathrm{tg}^{1\llbracket \mathbf{Pth}_{\boldsymbol{\mathcal{A}}^{(2)}}\rrbracket}_{s}\left(
\mathrm{tg}^{0\llbracket \mathbf{Pth}_{\boldsymbol{\mathcal{A}}^{(2)}}\rrbracket}_{s}\left(
\llbracket \mathfrak{P}^{(2)}\rrbracket_{s}
\right)\right)
&=
\mathrm{tg}^{0\llbracket \mathbf{Pth}_{\boldsymbol{\mathcal{A}}^{(2)}}\rrbracket}_{s}\left(
\llbracket \mathfrak{P}^{(2)}\rrbracket_{s}
\right);
\\
\mathrm{tg}^{0\llbracket \mathbf{Pth}_{\boldsymbol{\mathcal{A}}^{(2)}}\rrbracket}_{s}\left(
\mathrm{tg}^{1\llbracket \mathbf{Pth}_{\boldsymbol{\mathcal{A}}^{(2)}}\rrbracket}_{s}\left(
\llbracket \mathfrak{P}^{(2)}\rrbracket_{s}
\right)\right)
&=
\mathrm{tg}^{0\llbracket \mathbf{Pth}_{\boldsymbol{\mathcal{A}}^{(2)}}\rrbracket}_{s}\left(
\llbracket \mathfrak{P}^{(2)}\rrbracket_{s}
\right);
\\
\mathrm{tg}^{0\llbracket \mathbf{Pth}_{\boldsymbol{\mathcal{A}}^{(2)}}\rrbracket}_{s}\left(
\mathrm{sc}^{1\llbracket \mathbf{Pth}_{\boldsymbol{\mathcal{A}}^{(2)}}\rrbracket}_{s}\left(
\llbracket \mathfrak{P}^{(2)}\rrbracket_{s}
\right)\right)
&=
\mathrm{tg}^{0\llbracket \mathbf{Pth}_{\boldsymbol{\mathcal{A}}^{(2)}}\rrbracket}_{s}\left(
\llbracket \mathfrak{P}^{(2)}\rrbracket_{s}
\right).
\end{align*}
\end{proposition}
\begin{proof}
We will only present the proof for the first three equalities. The other ones are handled in a similar manner.

Regarding the first equality, the following chain of equalities holds.
\begin{flushleft}
$
\mathrm{sc}^{1\mathbf{Pth}_{\boldsymbol{\mathcal{A}}^{(2)}}}_{s}\left(
\mathrm{sc}^{0\mathbf{Pth}_{\boldsymbol{\mathcal{A}}^{(2)}}}_{s}\left(
\mathfrak{P}^{(2)}
\right)\right)
$
\allowdisplaybreaks
\begin{align*}
\qquad&=
\mathrm{sc}^{1\mathbf{Pth}_{\boldsymbol{\mathcal{A}}^{(2)}}}_{s}\left(
\mathrm{ip}^{(2,[1])\sharp}_{s}\left(
\mathrm{sc}^{0[\mathbf{PT}_{\boldsymbol{\mathcal{A}}}]}_{s}\left(
\mathrm{sc}^{([1],2)}_{s}\left(
\mathfrak{P}^{(2)}
\right)
\right)
\right)
\right)
\tag{1}
\\&=
\mathrm{ip}^{(2,[1])\sharp}_{s}\left(
\mathrm{sc}^{([1],2)}_{s}\left(
\mathrm{ip}^{(2,[1])\sharp}_{s}\left(
\mathrm{sc}^{0[\mathbf{PT}_{\boldsymbol{\mathcal{A}}}]}_{s}\left(
\mathrm{sc}^{([1],2)}_{s}\left(
\mathfrak{P}^{(2)}
\right)
\right)
\right)
\right)
\right)
\tag{2}
\\&=
\mathrm{ip}^{(2,[1])\sharp}_{s}\left(
\mathrm{sc}^{0[\mathbf{PT}_{\boldsymbol{\mathcal{A}}}]}_{s}\left(
\mathrm{sc}^{([1],2)}_{s}\left(
\mathfrak{P}^{(2)}
\right)
\right)
\right)
\tag{3}
\\&=
\mathrm{sc}^{0\mathbf{Pth}_{\boldsymbol{\mathcal{A}}^{(2)}}}_{s}\left(
\mathfrak{P}^{(2)}
\right).
\tag{4}
\end{align*}
\end{flushleft}

In the just stated chain of equalities, the first equality follows from Claim~\ref{CDPthCatAlgScZ}; the second equality recovers the description of the $1$-source operation symbol in the many-sorted partial $\Sigma^{\boldsymbol{\mathcal{A}}^{(2)}}$-algebra $\mathbf{Pth}_{\boldsymbol{\mathcal{A}}^{(2)}}$ according to Proposition~\ref{PDPthDCatAlg}; the third equality follows by Proposition~\ref{PDBasicEq}; finally, the last equality follows from Claim~\ref{CDPthCatAlgScZ}.

Note that the desired equality follows from the fact that the $\llbracket\cdot\rrbracket$-classes of the two second-order paths above are equal.

This completes the proof of the first equality.

Regarding the second equality, note that the following chain of equalities holds.
\begin{flushleft}
$
\mathrm{sc}^{0\mathbf{Pth}_{\boldsymbol{\mathcal{A}}^{(2)}}}_{s}\left(
\mathrm{sc}^{1\mathbf{Pth}_{\boldsymbol{\mathcal{A}}^{(2)}}}_{s}\left(
\mathfrak{P}^{(2)}
\right)
\right)
$
\allowdisplaybreaks
\begin{align*}
\qquad&=
\mathrm{sc}^{0\mathbf{Pth}_{\boldsymbol{\mathcal{A}}^{(2)}}}_{s}\left(
\mathrm{ip}^{(2,[1])\sharp}_{s}\left(
\mathrm{sc}^{([1],2)}_{s}\left(
\mathfrak{P}^{(2)}
\right)
\right)
\right)
\tag{1}
\\&=
\mathrm{ip}^{(2,[1])\sharp}_{s}\left(
\mathrm{sc}^{0[\mathbf{PT}_{\boldsymbol{\mathcal{A}}}]}_{s}\left(
\mathrm{sc}^{([1],2)}_{s}\left(
\mathrm{ip}^{(2,[1])\sharp}_{s}\left(
\mathrm{sc}^{([1],2)}_{s}\left(
\mathfrak{P}^{(2)}
\right)
\right)\right)
\right)\right)
\tag{2}
\\&=
\mathrm{ip}^{(2,[1])\sharp}_{s}\left(
\mathrm{sc}^{0[\mathbf{PT}_{\boldsymbol{\mathcal{A}}}]}_{s}\left(
\mathrm{sc}^{([1],2)}_{s}\left(
\mathfrak{P}^{(2)}
\right)
\right)\right)
\tag{3}
\\&=
\mathrm{sc}^{0\mathbf{Pth}_{\boldsymbol{\mathcal{A}}^{(2)}}}_{s}\left(
\mathfrak{P}^{(2)}
\right).
\tag{4}
\end{align*}
\end{flushleft}

In the just stated chain of equalities, the first equality recovers the description of the $1$-source operation symbol in the many-sorted partial $\Sigma^{\boldsymbol{\mathcal{A}}^{(2)}}$-algebra $\mathbf{Pth}_{\boldsymbol{\mathcal{A}}^{(2)}}$ according to Proposition~\ref{PDPthDCatAlg}; the second  equality follows from Claim~\ref{CDPthCatAlgScZ}; the third equality follows by Proposition~\ref{PDBasicEq}; finally, the last equality follows from Claim~\ref{CDPthCatAlgScZ}.

Note that the desired equality follows from the fact that the $\llbracket\cdot\rrbracket$-classes of the two second-order paths above are equal.

This completes the proof of the second equality.

Regarding the third equality, note that the following chain of equalities holds.
\begin{flushleft}
$
\mathrm{sc}^{0\mathbf{Pth}_{\boldsymbol{\mathcal{A}}^{(2)}}}_{s}\left(
\mathrm{tg}^{1\mathbf{Pth}_{\boldsymbol{\mathcal{A}}^{(2)}}}_{s}\left(
\mathfrak{P}^{(2)}
\right)
\right)
$
\allowdisplaybreaks
\begin{align*}
\qquad&=
\mathrm{sc}^{0\mathbf{Pth}_{\boldsymbol{\mathcal{A}}^{(2)}}}_{s}\left(
\mathrm{ip}^{(2,[1])\sharp}_{s}\left(
\mathrm{tg}^{([1],2)}_{s}\left(
\mathfrak{P}^{(2)}
\right)
\right)
\right)
\tag{1}
\\
&=
\mathrm{ip}^{(2,[1])\sharp}_{s}\left(
\mathrm{sc}^{0[\mathbf{PT}_{\boldsymbol{\mathcal{A}}}]}_{s}\left(
\mathrm{sc}^{([1],2)}_{s}\left(
\mathrm{ip}^{(2,[1])\sharp}_{s}\left(
\mathrm{tg}^{([1],2)}_{s}\left(
\mathfrak{P}^{(2)}
\right)
\right)\right)
\right)\right)
\tag{2}
\\
&=
\mathrm{ip}^{(2,[1])\sharp}_{s}\left(
\mathrm{sc}^{0[\mathbf{PT}_{\boldsymbol{\mathcal{A}}}]}_{s}\left(
\mathrm{tg}^{([1],2)}_{s}\left(
\mathfrak{P}^{(2)}
\right)
\right)\right)
\tag{3}
\\
&=
\mathrm{sc}^{0\mathbf{Pth}_{\boldsymbol{\mathcal{A}}^{(2)}}}_{s}\left(
\mathfrak{P}^{(2)}
\right).
\tag{4}
\end{align*}
\end{flushleft}

In the just stated chain of equalities, the first equality recovers the description of the $1$-target operation symbol in the many-sorted partial $\Sigma^{\boldsymbol{\mathcal{A}}^{(2)}}$-algebra $\mathbf{Pth}_{\boldsymbol{\mathcal{A}}^{(2)}}$ according to Proposition~\ref{PDPthDCatAlg}; the second  equality follows from Claim~\ref{CDPthCatAlgScZ}; the third equality follows by Proposition~\ref{PDBasicEq}; finally, the last equality follows from Claim~\ref{CDPthCatAlgScZ}.

Note that the desired equality follows from the fact that the $\llbracket\cdot\rrbracket$-classes of the two second-order paths above are equal.

This completes the proof of the third equality.

This completes the proof.
\end{proof}

\begin{proposition}\label{PDVVarAB2} Let $s$ be a sort in $S$ and $\llbracket \mathfrak{P}^{(2)} \rrbracket_{s}$, $\llbracket \mathfrak{Q}^{(2)} \rrbracket_{s}$ be  second-order path classes in $\llbracket \mathrm{Pth}_{\boldsymbol{\mathcal{A}}^{(2)}}\rrbracket_{s}$ satisfying that 
\[
\mathrm{sc}^{0\llbracket \mathbf{Pth}_{\boldsymbol{\mathcal{A}}^{(2)}}\rrbracket}_{s}\left(
\llbracket \mathfrak{Q}^{(2)}\rrbracket_{s}
\right)
=
\mathrm{tg}^{0\llbracket \mathbf{Pth}_{\boldsymbol{\mathcal{A}}^{(2)}}\rrbracket}_{s}\left(
\llbracket \mathfrak{P}^{(2)}\rrbracket_{s}
\right),
\]
then the following equalities hold
\allowdisplaybreaks
\begin{multline*}
\mathrm{sc}^{1\llbracket \mathbf{Pth}_{\boldsymbol{\mathcal{A}}^{(2)}}\rrbracket}_{s}\left(
\llbracket \mathfrak{Q}^{(2)}\rrbracket_{s}
\circ^{0\llbracket \mathbf{Pth}_{\boldsymbol{\mathcal{A}}^{(2)}}\rrbracket}_{s}
\llbracket \mathfrak{P}^{(2)}\rrbracket_{s}
\right)
\\=
\mathrm{sc}^{1\llbracket \mathbf{Pth}_{\boldsymbol{\mathcal{A}}^{(2)}}\rrbracket}_{s}\left(
\llbracket \mathfrak{Q}^{(2)}\rrbracket_{s}
\right)
\circ^{0\llbracket \mathbf{Pth}_{\boldsymbol{\mathcal{A}}^{(2)}}\rrbracket}_{s}
\mathrm{sc}^{1\llbracket \mathbf{Pth}_{\boldsymbol{\mathcal{A}}^{(2)}}\rrbracket}_{s}\left(
\llbracket \mathfrak{P}^{(2)}\rrbracket_{s}
\right);
\end{multline*}
\allowdisplaybreaks
\begin{multline*}
\mathrm{tg}^{1\llbracket \mathbf{Pth}_{\boldsymbol{\mathcal{A}}^{(2)}}\rrbracket}_{s}\left(
\llbracket \mathfrak{Q}^{(2)}\rrbracket_{s}
\circ^{0\llbracket \mathbf{Pth}_{\boldsymbol{\mathcal{A}}^{(2)}}\rrbracket}_{s}
\llbracket \mathfrak{P}^{(2)}\rrbracket_{s}
\right)
\\=
\mathrm{tg}^{1\llbracket \mathbf{Pth}_{\boldsymbol{\mathcal{A}}^{(2)}}\rrbracket}_{s}\left(
\llbracket \mathfrak{Q}^{(2)}\rrbracket_{s}
\right)
\circ^{0\llbracket \mathbf{Pth}_{\boldsymbol{\mathcal{A}}^{(2)}}\rrbracket}_{s}
\mathrm{tg}^{1\llbracket \mathbf{Pth}_{\boldsymbol{\mathcal{A}}^{(2)}}\rrbracket}_{s}\left(
\llbracket \mathfrak{P}^{(2)}\rrbracket_{s}
\right).
\end{multline*}
\end{proposition}
\begin{proof}
We will only present the proof for the first equality. The other one is handled in a similar manner.

Let us first prove that the following chain of equalities holds
\begin{flushleft}
$\mathrm{sc}^{1\mathbf{Pth}_{\boldsymbol{\mathcal{A}}^{(2)}}}_{s}\left(
\mathfrak{Q}^{(2)}
\circ^{0\mathbf{Pth}_{\boldsymbol{\mathcal{A}}^{(2)}}}_{s}
\mathfrak{P}^{(2)}
\right)
$
\allowdisplaybreaks
\begin{align*}
\qquad&=
\mathrm{ip}^{(2,[1])\sharp}_{s}\left(
\mathrm{sc}^{([1],2)}_{s}\left(
\mathfrak{Q}^{(2)}
\circ^{0\mathbf{Pth}_{\boldsymbol{\mathcal{A}}^{(2)}}}_{s}
\mathfrak{P}^{(2)}
\right)\right)
\tag{1}
\\&=
\mathrm{ip}^{(2,[1])\sharp}_{s}\left(
\mathrm{sc}^{([1],2)}_{s}\left(
\mathfrak{Q}^{(2)}
\right)
\circ^{0[\mathbf{PT}_{\boldsymbol{\mathcal{A}}}]}_{s}
\mathrm{sc}^{([1],2)}_{s}\left(
\mathfrak{P}^{(2)}
\right)\right)
\tag{2}
\\&=
\mathrm{ip}^{(2,[1])\sharp}_{s}\left(
\mathrm{sc}^{([1],2)}_{s}\left(
\mathfrak{Q}^{(2)}
\right)
\right)
\circ^{0\mathbf{Pth}_{\boldsymbol{\mathcal{A}}^{(2)}}}_{s}
\mathrm{ip}^{(2,[1])\sharp}_{s}\left(
\mathrm{sc}^{([1],2)}_{s}\left(
\mathfrak{P}^{(2)}
\right)\right)
\tag{3}
\\&=
\mathrm{sc}^{1\mathbf{Pth}_{\boldsymbol{\mathcal{A}}^{(2)}}}_{s}\left(
\mathfrak{Q}^{(2)}
\right)
\circ^{0\mathbf{Pth}_{\boldsymbol{\mathcal{A}}^{(2)}}}_{s}
\mathrm{sc}^{1\mathbf{Pth}_{\boldsymbol{\mathcal{A}}^{(2)}}}_{s}\left(
\mathfrak{P}^{(2)}
\right).
\tag{4}
\end{align*}
\end{flushleft}

In the just stated chain of equalities, the first equality unravels the description of the $1$-source operation symbol in the many-sorted partial $\Sigma^{\boldsymbol{\mathcal{A}}^{(2)}}$-algebra $\mathbf{Pth}_{\boldsymbol{\mathcal{A}}^{(2)}}$, according to Proposition~\ref{PDPthDCatAlg}; the second equality follows from Claim~\ref{CDPthCatAlgCompZ}; the third equality follows from Proposition~\ref{PDUIp}; finally, the last equality recovers the description of the $1$-source operation symbol in the many-sorted partial $\Sigma^{\boldsymbol{\mathcal{A}}^{(2)}}$-algebra $\mathbf{Pth}_{\boldsymbol{\mathcal{A}}^{(2)}}$, according to Proposition~\ref{PDPthDCatAlg}.

Note that the desired equality follows from the fact that the $\llbracket\cdot\rrbracket$-classes of the two second-order paths above are equal.

This completes the proof.
\end{proof}

\begin{proposition}\label{PDVVarAB3} Let $s$ be a sort in $S$ and $\llbracket \mathfrak{P}^{(2)} \rrbracket_{s}$, $\llbracket \mathfrak{P}'^{(2)} \rrbracket_{s}$, $\llbracket \mathfrak{Q}^{(2)} \rrbracket_{s}$, $\llbracket \mathfrak{Q}'^{(2)} \rrbracket_{s}$ be  second-order path classes in $\llbracket \mathrm{Pth}_{\boldsymbol{\mathcal{A}}^{(2)}}\rrbracket_{s}$ satisfying that 
\allowdisplaybreaks
\begin{align*}
\mathrm{sc}^{1\llbracket \mathbf{Pth}_{\boldsymbol{\mathcal{A}}^{(2)}}\rrbracket}_{s}\left(
\llbracket \mathfrak{Q}'^{(2)}\rrbracket_{s}
\right)
&=
\mathrm{tg}^{1\llbracket \mathbf{Pth}_{\boldsymbol{\mathcal{A}}^{(2)}}\rrbracket}_{s}\left(
\llbracket \mathfrak{Q}^{(2)}\rrbracket_{s}
\right);
\\
\mathrm{sc}^{1\llbracket \mathbf{Pth}_{\boldsymbol{\mathcal{A}}^{(2)}}\rrbracket}_{s}\left(
\llbracket \mathfrak{P}'^{(2)}\rrbracket_{s}
\right)
&=
\mathrm{tg}^{1\llbracket \mathbf{Pth}_{\boldsymbol{\mathcal{A}}^{(2)}}\rrbracket}_{s}\left(
\llbracket \mathfrak{P}^{(2)}\rrbracket_{s}
\right);
\\
\mathrm{sc}^{0\llbracket \mathbf{Pth}_{\boldsymbol{\mathcal{A}}^{(2)}}\rrbracket}_{s}\left(
\llbracket \mathfrak{Q}'^{(2)}\rrbracket_{s}
\right)
&=
\mathrm{tg}^{0\llbracket \mathbf{Pth}_{\boldsymbol{\mathcal{A}}^{(2)}}\rrbracket}_{s}\left(
\llbracket \mathfrak{P}'^{(2)}\rrbracket_{s}
\right);
\\
\mathrm{sc}^{0\llbracket \mathbf{Pth}_{\boldsymbol{\mathcal{A}}^{(2)}}\rrbracket}_{s}\left(
\llbracket \mathfrak{Q}^{(2)}\rrbracket_{s}
\right)
&=
\mathrm{tg}^{0\llbracket \mathbf{Pth}_{\boldsymbol{\mathcal{A}}^{(2)}}\rrbracket}_{s}\left(
\llbracket \mathfrak{P}^{(2)}\rrbracket_{s}
\right).
\end{align*}
then the following equality holds
\allowdisplaybreaks
\begin{multline*}
\left(
\llbracket \mathfrak{Q}'^{(2)}\rrbracket_{s}
\circ^{0\llbracket \mathbf{Pth}_{\boldsymbol{\mathcal{A}}^{(2)}}\rrbracket}_{s}
\llbracket \mathfrak{P}'^{(2)}\rrbracket_{s}
\right)
\circ^{1\llbracket \mathbf{Pth}_{\boldsymbol{\mathcal{A}}^{(2)}}\rrbracket}_{s}
\left(
\llbracket \mathfrak{Q}^{(2)}\rrbracket_{s}
\circ^{0\llbracket \mathbf{Pth}_{\boldsymbol{\mathcal{A}}^{(2)}}\rrbracket}_{s}
\llbracket \mathfrak{P}^{(2)}\rrbracket_{s}
\right)
\\=
\left(
\llbracket \mathfrak{Q}'^{(2)}\rrbracket_{s}
\circ^{1\llbracket \mathbf{Pth}_{\boldsymbol{\mathcal{A}}^{(2)}}\rrbracket}_{s}
\llbracket \mathfrak{Q}^{(2)}\rrbracket_{s}
\right)
\circ^{0\llbracket \mathbf{Pth}_{\boldsymbol{\mathcal{A}}^{(2)}}\rrbracket}_{s}
\left(
\llbracket \mathfrak{P}'^{(2)}\rrbracket_{s}
\circ^{1\llbracket \mathbf{Pth}_{\boldsymbol{\mathcal{A}}^{(2)}}\rrbracket}_{s}
\llbracket \mathfrak{P}^{(2)}\rrbracket_{s}
\right).
\end{multline*}
\end{proposition}
\begin{proof}
We will consider two cases according to the nature of the second-order paths $\mathfrak{Q}'^{(2)}$, $\mathfrak{Q}^{(2)}$, $\mathfrak{P}'^{(2)}$ and $\mathfrak{P}^{(2)}$. It could be the case that (1) all the second-order paths above are $(2,[1])$-identity second-order paths or (2) there exists at least one of them that is not a $(2,[1])$-identity second-order path.

For the first case, assume that, for suitable path terms $Q',Q,P',P$ in $\mathrm{PT}_{\boldsymbol{\mathcal{A}},s}$, it is the case that 
\allowdisplaybreaks
\begin{align*}
\mathfrak{Q}'^{(2)}&=
\mathrm{ip}^{(2,[1])\sharp}_{s}\left(
[Q']_{s}
\right);
&
\mathfrak{P}'^{(2)}&=
\mathrm{ip}^{(2,[1])\sharp}_{s}\left(
[P']_{s}
\right);
\\
\mathfrak{Q}^{(2)}&=
\mathrm{ip}^{(2,[1])\sharp}_{s}\left(
[Q]_{s}
\right);
&
\mathfrak{P}^{(2)}&=
\mathrm{ip}^{(2,[1])\sharp}_{s}\left(
[P]_{s}
\right).
\end{align*}

Note that, in order to satisfy the first equality that holds by hypothesis, the following chain of equalities must hold
\begin{align*}
\Bigl\llbracket \mathrm{ip}^{(2,[1])\sharp}_{s}\left(
\mathrm{sc}^{([1],2)}_{s}\left(
\mathfrak{Q}'^{(2)}
\right)
\right) \Bigr\rrbracket_{s}
&=
\Bigl\llbracket 
\mathrm{sc}^{1\mathbf{Pth}_{\boldsymbol{\mathcal{A}}^{(2)}}}_{s}\left(
\mathfrak{Q}'^{(2)}
\right)
\Bigr\rrbracket_{s}
\tag{1}
\\
&=
\mathrm{sc}^{1\llbracket \mathbf{Pth}_{\boldsymbol{\mathcal{A}}^{(2)}} \rrbracket}_{s}\left(
\llbracket
\mathfrak{Q}'^{(2)}
\rrbracket
\right)
\tag{2}
\\&=
\mathrm{tg}^{1\llbracket \mathbf{Pth}_{\boldsymbol{\mathcal{A}}^{(2)}} \rrbracket}_{s}\left(
\llbracket
\mathfrak{Q}^{(2)}
\rrbracket
\right)
\tag{3}
\\&=
\Bigl\llbracket 
\mathrm{tg}^{1\mathbf{Pth}_{\boldsymbol{\mathcal{A}}^{(2)}}}_{s}\left(
\mathfrak{Q}^{(2)}
\right)
\Bigr\rrbracket_{s}
\tag{4}
\\&=
\Bigl\llbracket \mathrm{ip}^{(2,[1])\sharp}_{s}\left(
\mathrm{tg}^{([1],2)}_{s}\left(
\mathfrak{Q}^{(2)}
\right)
\right) \Bigr\rrbracket_{s}.
\tag{5}
\end{align*}

In the just stated chain of equalities, the first equality recovers the description of the $1$-source operation symbol in the many-sorted partial $\Sigma^{\boldsymbol{\mathcal{A}}^{(2)}}$-algebra $\mathbf{Pth}_{\boldsymbol{\mathcal{A}}^{(2)}}$, according to Proposition~\ref{PDPthDCatAlg}; the second equality recovers the description of the $1$-source operation symbol in the many-sorted partial $\Sigma^{\boldsymbol{\mathcal{A}}^{(2)}}$-algebra $\llbracket \mathbf{Pth}_{\boldsymbol{\mathcal{A}}^{(2)}}\rrbracket $, according to Proposition~\ref{PDVDCatAlg}; the third equality holds by assumption; the fourth equality unravels the description of the $1$-target operation symbol in the many-sorted partial $\Sigma^{\boldsymbol{\mathcal{A}}^{(2)}}$-algebra $\llbracket \mathbf{Pth}_{\boldsymbol{\mathcal{A}}^{(2)}}\rrbracket $, according to Proposition~\ref{PDVDCatAlg}; finally, the last equality recovers the description of the $1$-target operation symbol in the many-sorted partial $\Sigma^{\boldsymbol{\mathcal{A}}^{(2)}}$-algebra $\mathbf{Pth}_{\boldsymbol{\mathcal{A}}^{(2)}}$, according to Proposition~\ref{PDPthDCatAlg}.

Taking into account the above equality and Proposition~\ref{PDVDUId}, we conclude that 
\[
\mathrm{ip}^{(2,[1])\sharp}_{s}\left(
\mathrm{sc}^{([1],2)}_{s}\left(
\mathfrak{Q}'^{(2)}
\right)
\right)
=
\mathrm{ip}^{(2,[1])\sharp}_{s}\left(
\mathrm{tg}^{([1],2)}_{s}\left(
\mathfrak{Q}^{(2)}
\right)
\right).
\]

Applying the $([1],2)$-source mapping on both sides above we conclude, according to Proposition~\ref{PDBasicEq}, that 
\[
[Q']_{s}=
\mathrm{sc}^{([1],2)}_{s}\left(
\mathfrak{Q}'^{(2)}
\right)
=
\mathrm{tg}^{([1],2)}_{s}\left(
\mathfrak{Q}^{(2)}
\right)
=
[Q]_{s}.
\]

Consequently, the two second-order paths $\mathfrak{Q}^{(2)}$ and $\mathfrak{Q}'^{(2)}$ are equal, since
\[
\mathfrak{Q}^{(2)}
=\mathrm{ip}^{(2,[1])\sharp}_{s}\left([Q']_{s}\right)
=\mathrm{ip}^{(2,[1])\sharp}_{s}\left([Q]_{s}\right)
=\mathfrak{Q}^{(2)}.
\]

Taking into account that, by assumption, the equality 
\[
\mathrm{sc}^{1\llbracket\mathbf{Pth}_{\boldsymbol{\mathcal{A}}^{(2)}} \rrbracket}_{s}
\left(
\llbracket 
\mathfrak{P}'^{(2)}
\rrbracket 
\right)
=
\mathrm{tg}^{1\llbracket\mathbf{Pth}_{\boldsymbol{\mathcal{A}}^{(2)}} \rrbracket}_{s}
\left(
\llbracket 
\mathfrak{P}^{(2)}
\rrbracket 
\right)
\]
also holds, we conclude, by a similar argument to that above, that 
$[P']_{s}=[P]_{s}$ and 
$\mathfrak{P}'^{(2)}=\mathfrak{P}^{(2)}$.

Then the following chain of equalities holds
\begin{flushleft}
$
\left(
\mathfrak{Q}'^{(2)}
\circ^{0\mathbf{Pth}_{\boldsymbol{\mathcal{A}}^{(2)}}}_{s}
\mathfrak{P}'^{(2)}
\right)
\circ^{1\mathbf{Pth}_{\boldsymbol{\mathcal{A}}^{(2)}}}_{s}
\left(
\mathfrak{Q}^{(2)}
\circ^{0\mathbf{Pth}_{\boldsymbol{\mathcal{A}}^{(2)}}}_{s}
\mathfrak{P}^{(2)}
\right)
$
\allowdisplaybreaks
\begin{align*}
&=
\left(
\mathrm{ip}^{(2,[1])\sharp}_{s}\left(
[Q]_{s}
\right)
\circ^{0\mathbf{Pth}_{\boldsymbol{\mathcal{A}}^{(2)}}}_{s}
\mathrm{ip}^{(2,[1])\sharp}_{s}\left(
[P]_{s}
\right)
\right)
\circ^{1\mathbf{Pth}_{\boldsymbol{\mathcal{A}}^{(2)}}}_{s}
\\&\qquad\qquad\qquad\qquad\qquad
\left(
\mathrm{ip}^{(2,[1])\sharp}_{s}\left(
[Q]_{s}
\right)
\circ^{0\mathbf{Pth}_{\boldsymbol{\mathcal{A}}^{(2)}}}_{s}
\mathrm{ip}^{(2,[1])\sharp}_{s}\left(
[P]_{s}
\right)
\right)
\tag{1}
\\&=
\mathrm{ip}^{(2,[1])\sharp}_{s}\left(
\left[
Q\circ^{0\mathbf{PT}_{\boldsymbol{\mathcal{A}}} }_{s}
P
\right]_{s}
\right)
\circ^{1\mathbf{Pth}_{\boldsymbol{\mathcal{A}}^{(2)}}}_{s}
\mathrm{ip}^{(2,[1])\sharp}_{s}\left(
\left[
Q\circ^{0\mathbf{PT}_{\boldsymbol{\mathcal{A}}} }_{s}
P
\right]_{s}
\right)
\tag{2}
\\&=
\mathrm{ip}^{(2,[1])\sharp}_{s}\left(
\left[
Q\circ^{0\mathbf{PT}_{\boldsymbol{\mathcal{A}}} }_{s}
P
\right]_{s}
\right)
\tag{3}
\\&=
\mathrm{ip}^{(2,[1])\sharp}_{s}\left(
[Q]_{s}
\right)
\circ^{0\mathbf{Pth}_{\boldsymbol{\mathcal{A}}^{(2)}}}_{s}
\mathrm{ip}^{(2,[1])\sharp}_{s}\left(
[P]_{s}
\right)
\tag{4}
\\&=
\mathfrak{Q}^{(2)}
\circ^{0\mathbf{Pth}_{\boldsymbol{\mathcal{A}}^{(2)}}}_{s}
\mathfrak{P}^{(2)}
\tag{5}
\\&=
\left(
\mathfrak{Q}'^{(2)}
\circ^{1\mathbf{Pth}_{\boldsymbol{\mathcal{A}}^{(2)}}}_{s}
\mathfrak{Q}^{(2)}
\right)
\circ^{0\mathbf{Pth}_{\boldsymbol{\mathcal{A}}^{(2)}}}_{s}
\left(
\mathfrak{P}'^{(2)}
\circ^{1\mathbf{Pth}_{\boldsymbol{\mathcal{A}}^{(2)}}}_{s}
\mathfrak{P}^{(2)}
\right).
\tag{6}
\end{align*}
\end{flushleft}

In the just stated chain of equalities, the first equality follows from the description of the second-order paths as $(2,[1])$-identity second-order paths; the second equality follows from Proposition~\ref{PDUIp}; the third equality follows from the fact that, according to Definition~\ref{DDPthComp}, the $(2,[1])$-identity second-order paths are idempotent for the $1$-composition; the fourth equality follows from Proposition~\ref{PDUIp}; the fifth equality recovers the description of the second-order paths $\mathfrak{P}^{(2)}$ and $\mathfrak{Q}^{(2)}$; finally, the last equality follows from the fact that, according to Definition~\ref{DDPthComp}, the $(2,[1])$-identity second-order paths are idempotent for the $1$-composition. In this regard, let us recall that, by the previous discussion, $\mathfrak{Q}'^{(2)}=\mathfrak{Q}^{(2)}$ and $\mathfrak{P}'^{(2)}=\mathfrak{P}^{(2)}$.

Note that the desired equality follows from the fact that the $\llbracket\cdot\rrbracket$-classes of the two second-order paths above are equal.

This completes case~(1), where all second-order paths above are $(2,[1])$-identity second-order paths.

Now, consider the case~(2), i.e., the case in which at least one of the second-order paths under consideration is not a $(2,[1])$-identity second-order path.

Consider, on one hand, the second-order path
\[
\left(
\mathfrak{Q}'^{(2)}
\circ^{0\mathbf{Pth}_{\boldsymbol{\mathcal{A}}^{(2)}}}_{s}
\mathfrak{P}'^{(2)}
\right)
\circ^{1\mathbf{Pth}_{\boldsymbol{\mathcal{A}}^{(2)}}}_{s}
\left(
\mathfrak{Q}^{(2)}
\circ^{0\mathbf{Pth}_{\boldsymbol{\mathcal{A}}^{(2)}}}_{s}
\mathfrak{P}^{(2)}
\right).
\]

An schematic representation of this second-order path is included in Figure~\ref{FDVVarAB3Lft}. Note that we have removed the superscripts $[\mathbf{PT}_{\boldsymbol{\mathcal{A}}}]$ for clarity.  In virtue of Corollary~\ref{CDPthWB}, this is a coherent head-constant echelonless second-order path associated to the $0$-composition operation symbol. Note that the second-order path extraction procedure from Lemma~\ref{LDPthExtract} applied to it, retrieves the second-order paths $\mathfrak{Q}'^{(2)}\circ^{1\mathbf{Pth}_{\boldsymbol{\mathcal{A}}^{(2)}}} \mathfrak{Q}^{(2)}$ and $\mathfrak{P}'^{(2)}\circ^{1\mathbf{Pth}_{\boldsymbol{\mathcal{A}}^{(2)}}} \mathfrak{P}^{(2)}$, respectively. 

Therefore, the second-order Curry-Howard mapping at this second-order path is given by the following expression.
\allowdisplaybreaks
\begin{multline*}
\mathrm{CH}^{(2)}_{s}\left(
\left(
\mathfrak{Q}'^{(2)}
\circ^{0\mathbf{Pth}_{\boldsymbol{\mathcal{A}}^{(2)}}}_{s}
\mathfrak{P}'^{(2)}
\right)
\circ^{1\mathbf{Pth}_{\boldsymbol{\mathcal{A}}^{(2)}}}_{s}
\left(
\mathfrak{Q}^{(2)}
\circ^{0\mathbf{Pth}_{\boldsymbol{\mathcal{A}}^{(2)}}}_{s}
\mathfrak{P}^{(2)}
\right)
\right)
\\=
\mathrm{CH}^{(2)}_{s}\left(
\mathfrak{Q}'^{(2)}\circ^{1\mathbf{Pth}_{\boldsymbol{\mathcal{A}}^{(2)}}} \mathfrak{Q}^{(2)}
\right)
\circ^{0\mathbf{T}_{\Sigma^{\boldsymbol{\mathcal{A}}^{(2)}}}(X)}_{s}
\mathrm{CH}^{(2)}_{s}\left(
\mathfrak{P}'^{(2)}\circ^{1\mathbf{Pth}_{\boldsymbol{\mathcal{A}}^{(2)}}} \mathfrak{P}^{(2)}
\right).
\end{multline*}

\begin{figure}
\centering
\begin{tikzpicture}[
fletxa/.style={thick, ->}]
\tikzstyle{every state}=[very thick, draw=blue!50,fill=blue!20, inner sep=2pt,minimum size=20pt]

\node[] (a1) at (0,0) [] {$\mathrm{sc}^{([1],2)}_{s}\left(\mathfrak{Q}^{(2)}\right)$}; 
\node[] (a2) at (1.6,0) []  {$\circ^{0}_{s}$};
\node[] (a3) at (3.2,0) [] {$\mathrm{sc}^{([1],2)}_{s}\left(\mathfrak{P}^{(2)}\right)$};

\node[] (b1) at (0,-1) [] {\rotatebox{-90}{$\,=$}}; 
\node[] (b2) at (1.6,-1) [] {$\vdots$};
\node[] (b3) at (3.5,-.95) [] {\rotatebox{-90}{
$\xymatrix@C=30pt{
\,
\ar@{=>}[r]^-{\text{
\rotatebox{90}{$\mathfrak{P}^{(2)}$}
}}
&
\,
}$}}; 

\node[] (c1) at (0,-2) [] {$\mathrm{sc}^{([1],2)}_{s}\left(\mathfrak{Q}^{(2)}\right)$}; 
\node[] (c2) at (1.6,-2) []   {$\circ^{0}_{s}$};
\node[] (c3) at (3.2,-2) [] {$\mathrm{tg}^{([1],2)}_{s}\left(\mathfrak{P}^{(2)}\right)$}; 

\node[] (d1) at (0.3,-2.95) [] {\rotatebox{-90}{
$\xymatrix@C=30pt{
\,
\ar@{=>}[r]^-{\text{
\rotatebox{90}{$\mathfrak{Q}^{(2)}$}
}}
&
\,
}$}}; 
\node[] (d2) at (1.6,-3) [] {$\vdots$};
\node[] (d3) at (3.2,-2.95) []  {\rotatebox{-90}{$\,=$}};

\node[] (e1) at (0,-4) [] {$\mathrm{tg}^{([1],2)}_{s}\left(\mathfrak{Q}^{(2)}\right)$}; 
\node[] (e2) at (1.6,-4) []   {$\circ^{0}_{s}$};
\node[] (e3) at (3.2,-4) [] {$\mathrm{sc}^{([1],2)}_{s}\left(\mathfrak{P}'^{(2)}\right)$}; 

\node[] (f1) at (0,-5) [] {\rotatebox{-90}{$\,=$}}; 
\node[] (f2) at (1.6,-5) [] {$\vdots$};
\node[] (f3) at (3.5,-4.95) []  {\rotatebox{-90}{
$\xymatrix@C=30pt{
\,
\ar@{=>}[r]^-{\text{
\rotatebox{90}{$\mathfrak{P}'^{(2)}$}
}}
&
\,
}$}}; 

\node[] (g1) at (0,-6) [] {$\mathrm{sc}^{([1],2)}_{s}\left(\mathfrak{Q}'^{(2)}\right)$}; 
\node[] (g2) at (1.6,-6) []   {$\circ^{0}_{s}$};
\node[] (g3) at (3.2,-6) [] {$\mathrm{tg}^{([1],2)}_{s}\left(\mathfrak{P}'^{(2)}\right)$}; 

\node[] (h1) at (0.3,-6.95) [] 
{\rotatebox{-90}{
$\xymatrix@C=30pt{
\,
\ar@{=>}[r]^-{\text{
\rotatebox{90}{$\mathfrak{Q}'^{(2)}$}
}}
&
\,
}$}}; 
\node[] (h2) at (1.6,-7) [] {$\vdots$};
\node[] (h3) at (3.2,-6.95) [] {\rotatebox{-90}{$\,=$}};

\node[] (i1) at (0,-8) [] {$\mathrm{tg}^{([1],2)}_{s}\left(\mathfrak{Q}'^{(2)}\right)$}; 
\node[] (i2) at (1.6,-8) []  {$\circ^{0}_{s}$};
\node[] (i3) at (3.2,-8) [] {$\mathrm{tg}^{([1],2)}_{s}\left(\mathfrak{P}'^{(2)}\right)$}; 

\tikzset{encercla/.style={draw=black, line width=.5pt, inner sep=0pt, rectangle, rounded corners}};
\node [encercla,  fit=(a1)(c1) ] {} ;
\node [encercla,  fit=(e1)(g1) ] {} ;
\node [encercla,  fit=(i1) ] {} ;
\node [encercla,  fit=(a3) ] {} ;
\node [encercla,  fit=(c3)(e3) ] {} ;
\node [encercla,  fit=(g3)(i3) ] {} ;
\end{tikzpicture}
\caption{
$(
\mathfrak{Q}'^{(2)}
\circ^{0\mathbf{Pth}_{\boldsymbol{\mathcal{A}}^{(2)}}}_{s}
\mathfrak{P}'^{(2)}
)
\circ^{1\mathbf{Pth}_{\boldsymbol{\mathcal{A}}^{(2)}}}_{s}
(
\mathfrak{Q}^{(2)}
\circ^{0\mathbf{Pth}_{\boldsymbol{\mathcal{A}}^{(2)}}}_{s}
\mathfrak{P}^{(2)}
)$.
}\label{FDVVarAB3Lft}
\end{figure}

Now consider, on the other hand, the second-order path 
\[
\left(
\mathfrak{Q}'^{(2)}
\circ^{1\mathbf{Pth}_{\boldsymbol{\mathcal{A}}^{(2)}}}_{s}
\mathfrak{Q}^{(2)}
\right)
\circ^{0\mathbf{Pth}_{\boldsymbol{\mathcal{A}}^{(2)}}}_{s}
\left(
\mathfrak{P}'^{(2)}
\circ^{1\mathbf{Pth}_{\boldsymbol{\mathcal{A}}^{(2)}}}_{s}
\mathfrak{P}^{(2)}
\right).
\]

An schematic representation of this second-order path is included in Figure~\ref{FDVVarAB3Rgt}. Note that, as before, we have removed the superscripts $[\mathbf{PT}_{\boldsymbol{\mathcal{A}}}]$ for clarity.  This is, in virtue of Corollary~\ref{CDPthWB}, a coherent head-constant echelonless second-order path associated to the $0$-composition operation symbol. Note that the second-order path extraction procedure from Lemma~\ref{LDPthExtract} applied to it, retrieves the second-order paths $\mathfrak{Q}'^{(2)}\circ^{1\mathbf{Pth}_{\boldsymbol{\mathcal{A}}^{(2)}}} \mathfrak{Q}^{(2)}$ and $\mathfrak{P}'^{(2)}\circ^{1\mathbf{Pth}_{\boldsymbol{\mathcal{A}}^{(2)}}} \mathfrak{P}^{(2)}$, respectively. 

Therefore, the second-order Curry-Howard mapping at this second-order path is given by the following expression.
\allowdisplaybreaks
\begin{multline*}
\mathrm{CH}^{(2)}_{s}\left(
\left(
\mathfrak{Q}'^{(2)}
\circ^{1\mathbf{Pth}_{\boldsymbol{\mathcal{A}}^{(2)}}}_{s}
\mathfrak{Q}^{(2)}
\right)
\circ^{0\mathbf{Pth}_{\boldsymbol{\mathcal{A}}^{(2)}}}_{s}
\left(
\mathfrak{P}'^{(2)}
\circ^{1\mathbf{Pth}_{\boldsymbol{\mathcal{A}}^{(2)}}}_{s}
\mathfrak{P}^{(2)}
\right)
\right)
\\=
\mathrm{CH}^{(2)}_{s}\left(
\mathfrak{Q}'^{(2)}\circ^{1\mathbf{Pth}_{\boldsymbol{\mathcal{A}}^{(2)}}} \mathfrak{Q}^{(2)}
\right)
\circ^{0\mathbf{T}_{\Sigma^{\boldsymbol{\mathcal{A}}^{(2)}}}(X)}_{s}
\mathrm{CH}^{(2)}_{s}\left(
\mathfrak{P}'^{(2)}\circ^{1\mathbf{Pth}_{\boldsymbol{\mathcal{A}}^{(2)}}} \mathfrak{P}^{(2)}
\right).
\end{multline*}

\begin{figure}
\centering
\begin{tikzpicture}[
fletxa/.style={thick, ->}]
\tikzstyle{every state}=[very thick, draw=blue!50,fill=blue!20, inner sep=2pt,minimum size=20pt]

\node[] (a1) at (0,0) [] {$\mathrm{sc}^{([1],2)}_{s}\left(\mathfrak{Q}^{(2)}\right)$}; 
\node[] (a2) at (1.6,0) []  {$\circ^{0}_{s}$};
\node[] (a3) at (3.2,0) [] {$\mathrm{sc}^{([1],2)}_{s}\left(\mathfrak{P}^{(2)}\right)$};

\node[] (b1) at (0,-1) [] {\rotatebox{-90}{$\,=$}}; 
\node[] (b2) at (1.6,-1) [] {$\vdots$};
\node[] (b3) at (3.5,-.95) [] {\rotatebox{-90}{
$\xymatrix@C=30pt{
\,
\ar@{=>}[r]^-{\text{
\rotatebox{90}{$\mathfrak{P}^{(2)}$}
}}
&
\,
}$}}; 

\node[] (c1) at (0,-2) [] {$\vdots$}; 
\node[] (c2) at (1.6,-2) []  {$\circ^{0}_{s}$};
\node[] (c3) at (3.2,-2) [] {$\mathrm{tg}^{([1],2)}_{s}\left(\mathfrak{P}^{(2)}\right)$}; 

\node[] (d1) at (0,-3) [] {\rotatebox{-90}{$\,=$}};  
\node[] (d2) at (1.6,-3) [] {$\vdots$};
\node[] (d3) at (3.2,-3) [] {\rotatebox{-90}{$\,=$}}; 

\node[] (e1) at (0,-4) [] {$\vdots$}; 
\node[] (e2) at (1.6,-4) []  {$\circ^{0}_{s}$};
\node[] (e3) at (3.2,-4) [] {$\mathrm{sc}^{([1],2)}_{s}\left(\mathfrak{P}'^{(2)}\right)$};

\node[] (f1) at (0,-5) [] {\rotatebox{-90}{$\,=$}}; 
\node[] (f2) at (1.6,-5) [] {$\vdots$};
\node[] (f3) at (3.5,-4.95) []  {\rotatebox{-90}{
$\xymatrix@C=30pt{
\,
\ar@{=>}[r]^-{\text{
\rotatebox{90}{$\mathfrak{P}'^{(2)}$}
}}
&
\,
}$}};

\node[] (g1) at (0,-6) [] {$\mathrm{sc}^{([1],2)}_{s}\left(\mathfrak{Q}^{(2)}\right)$}; 
\node[] (g2) at (1.6,-6) []  {$\circ^{0}_{s}$};
\node[] (g3) at (3.2,-6) [] {$\mathrm{tg}^{([1],2)}_{s}\left(\mathfrak{P}'^{(2)}\right)$}; 

\node[] (h1) at (.3,-6.95) [] {\rotatebox{-90}{
$\xymatrix@C=30pt{
\,
\ar@{=>}[r]^-{\text{
\rotatebox{90}{$\mathfrak{Q}^{(2)}$}
}}
&
\,
}$}}; 
\node[] (h2) at (1.6,-7) [] {$\vdots$};
\node[] (h3) at (3.2,-7) [] {\rotatebox{-90}{$\,=$}};

\node[] (i1) at (0,-8) [] {$\mathrm{tg}^{([1],2)}_{s}\left(\mathfrak{Q}^{(2)}\right)$};
\node[] (i2) at (1.6,-8) []{$\circ^{0}_{s}$};
\node[] (i3) at (3.2,-8) [] {$\vdots$}; 

\node[] (j1) at (0,-9) [] {\rotatebox{-90}{$\,=$}};  
\node[] (j2) at (1.6,-9) [] {$\vdots$};
\node[] (j3) at (3.2,-9) [] {\rotatebox{-90}{$\,=$}}; 

\node[] (k1) at (0,-10) [] {$\mathrm{sc}^{([1],2)}_{s}\left(\mathfrak{Q}'^{(2)}\right)$}; 
\node[] (k2) at (1.6,-10) [] {$\circ^{0}_{s}$};
\node[] (k3) at (3.2,-10) [] {$\vdots$}; 

\node[] (l1) at (0.3,-10.95) []  {\rotatebox{-90}{
$\xymatrix@C=30pt{
\,
\ar@{=>}[r]^-{\text{
\rotatebox{90}{$\mathfrak{P}'^{(2)}$}
}}
&
\,
}$}}; 
\node[] (l2) at (1.6,-11) [] {$\vdots$};
\node[] (l3) at (3.2,-11) [] {\rotatebox{-90}{$\,=$}};

\node[] (m1) at (0,-12) [] {$\mathrm{tg}^{([1],2)}_{s}\left(\mathfrak{Q}'^{(2)}\right)$}; 
\node[] (m2) at (1.6,-12) []  {$\circ^{0}_{s}$};
\node[] (m3) at (3.2,-12) [] {$\mathrm{tg}^{([1],2)}_{s}\left(\mathfrak{P}'^{(2)}\right)$}; 

\tikzset{encercla/.style={draw=black, line width=.5pt, inner sep=0pt, rectangle, rounded corners}};
\node [encercla,  fit=(a1)(g1) ] {} ;
\node [encercla,  fit=(i1)(k1) ] {} ;
\node [encercla,  fit=(m1) ] {} ;
\node [encercla,  fit=(a3) ] {} ;
\node [encercla,  fit=(c3)(e3) ] {} ;
\node [encercla,  fit=(g3)(m3) ] {} ;
\end{tikzpicture}
\caption{
$(
\mathfrak{Q}'^{(2)}
\circ^{1\mathbf{Pth}_{\boldsymbol{\mathcal{A}}^{(2)}}}_{s}
\mathfrak{Q}^{(2)}
)
\circ^{0\mathbf{Pth}_{\boldsymbol{\mathcal{A}}^{(2)}}}_{s}
(
\mathfrak{P}'^{(2)}
\circ^{1\mathbf{Pth}_{\boldsymbol{\mathcal{A}}^{(2)}}}_{s}
\mathfrak{P}^{(2)}
)$.
}\label{FDVVarAB3Rgt}
\end{figure}

All in all, we conclude that the above second-order paths are in the kernel of the second-order Curry-Howard mapping, i.e.,
\allowdisplaybreaks
\begin{multline*}
\left[
\left(
\mathfrak{Q}'^{(2)}
\circ^{0\mathbf{Pth}_{\boldsymbol{\mathcal{A}}^{(2)}}}_{s}
\mathfrak{P}'^{(2)}
\right)
\circ^{1\mathbf{Pth}_{\boldsymbol{\mathcal{A}}^{(2)}}}_{s}
\left(
\mathfrak{Q}^{(2)}
\circ^{0\mathbf{Pth}_{\boldsymbol{\mathcal{A}}^{(2)}}}_{s}
\mathfrak{P}^{(2)}
\right)
\right]_{s}
\\=
\left[
\left(
\mathfrak{Q}'^{(2)}
\circ^{1\mathbf{Pth}_{\boldsymbol{\mathcal{A}}^{(2)}}}_{s}
\mathfrak{Q}^{(2)}
\right)
\circ^{0\mathbf{Pth}_{\boldsymbol{\mathcal{A}}^{(2)}}}_{s}
\left(
\mathfrak{P}'^{(2)}
\circ^{1\mathbf{Pth}_{\boldsymbol{\mathcal{A}}^{(2)}}}_{s}
\mathfrak{P}^{(2)}
\right)
\right]_{s}.
\end{multline*}

Therefore, their $\llbracket\cdot \rrbracket$-class must coincide, i.e.,
\allowdisplaybreaks
\begin{multline*}
\biggl\llbracket
\left(
\mathfrak{Q}'^{(2)}
\circ^{0\mathbf{Pth}_{\boldsymbol{\mathcal{A}}^{(2)}}}_{s}
\mathfrak{P}'^{(2)}
\right)
\circ^{1\mathbf{Pth}_{\boldsymbol{\mathcal{A}}^{(2)}}}_{s}
\left(
\mathfrak{Q}^{(2)}
\circ^{0\mathbf{Pth}_{\boldsymbol{\mathcal{A}}^{(2)}}}_{s}
\mathfrak{P}^{(2)}
\right)
\biggr\rrbracket_{s}
\\=
\biggl\llbracket
\left(
\mathfrak{Q}'^{(2)}
\circ^{1\mathbf{Pth}_{\boldsymbol{\mathcal{A}}^{(2)}}}_{s}
\mathfrak{Q}^{(2)}
\right)
\circ^{0\mathbf{Pth}_{\boldsymbol{\mathcal{A}}^{(2)}}}_{s}
\left(
\mathfrak{P}'^{(2)}
\circ^{1\mathbf{Pth}_{\boldsymbol{\mathcal{A}}^{(2)}}}_{s}
\mathfrak{P}^{(2)}
\right)
\biggr\rrbracket_{s}.
\end{multline*}

This completes the proof.
\end{proof}

In virtue of the foregoing propositions, we can consider the $S$-sorted category given by second-order path $\llbracket \cdot \rrbracket$-classes.

\begin{restatable}{definition}{DDVCat}
\label{DDVCat} 
\index{path!second-order!$\llbracket \mathsf{Pth}_{\boldsymbol{\mathcal{A}}^{(2)}}\rrbracket$}
Let $\llbracket \mathsf{Pth}_{\boldsymbol{\mathcal{A}}^{(2)}}\rrbracket$ denote the ordered tuple
\allowdisplaybreaks
\begin{multline*}
\llbracket \mathsf{Pth}_{\boldsymbol{\mathcal{A}}^{(2)}}\rrbracket=
\left(
\llbracket \mathbf{Pth}_{\boldsymbol{\mathcal{A}}^{(2)}}\rrbracket, 
\left(
\circ^{0 \llbracket \mathbf{Pth}_{\boldsymbol{\mathcal{A}}^{(2)}}\rrbracket},
\mathrm{sc}^{0\llbracket \mathbf{Pth}_{\boldsymbol{\mathcal{A}}^{(2)}}\rrbracket},
\mathrm{tg}^{0\llbracket \mathbf{Pth}_{\boldsymbol{\mathcal{A}}^{(2)}}\rrbracket},
\right),
\right.
\\
\left.
\left(
\circ^{1 \llbracket \mathbf{Pth}_{\boldsymbol{\mathcal{A}}^{(2)}}\rrbracket},
\mathrm{sc}^{1\llbracket \mathbf{Pth}_{\boldsymbol{\mathcal{A}}^{(2)}}\rrbracket},
\mathrm{tg}^{1\llbracket \mathbf{Pth}_{\boldsymbol{\mathcal{A}}^{(2)}}\rrbracket},
\right)
\right).
\end{multline*}
\end{restatable}

\begin{restatable}{proposition}{PDVCat}
\label{PDVCat} $\llbracket \mathsf{Pth}_{\boldsymbol{\mathcal{A}}^{(2)}}\rrbracket$ is an $S$-sorted $2$-category.
\end{restatable}
\begin{proof}
Let us recall that, for every sort $s\in S$, the structure $(\llbracket \mathbf{Pth}_{\boldsymbol{\mathcal{A}}^{(2)},s}\rrbracket, \xi_{0, s}, \xi_{1, s})$, where 
\begin{align*}
\xi_{0,s}&=\left(
\circ^{0 \llbracket \mathbf{Pth}_{\boldsymbol{\mathcal{A}}^{(2)}}\rrbracket}_{s},
\mathrm{sc}^{0\llbracket \mathbf{Pth}_{\boldsymbol{\mathcal{A}}^{(2)}}\rrbracket}_{s},
\mathrm{tg}^{0\llbracket \mathbf{Pth}_{\boldsymbol{\mathcal{A}}^{(2)}}\rrbracket}_{s},
\right)
\\
\xi_{1,s}&=
\left(
\circ^{1 \llbracket \mathbf{Pth}_{\boldsymbol{\mathcal{A}}^{(2)}}\rrbracket}_{s},
\mathrm{sc}^{1\llbracket \mathbf{Pth}_{\boldsymbol{\mathcal{A}}^{(2)}}\rrbracket}_{s},
\mathrm{tg}^{1\llbracket \mathbf{Pth}_{\boldsymbol{\mathcal{A}}^{(2)}}\rrbracket}_{s},
\right)
\end{align*}
is a single-sorted  $2$-category in virtue of Definition~\ref{D2Cat} and Propositions~\ref{PDVVarA2}, \ref{PDVVarA3}, \ref{PDVVarA4}, \ref{PDVVarA5}, \ref{PDVVarA6}, \ref{PDVVarB2}, \ref{PDVVarB3}, \ref{PDVVarB4}, \ref{PDVVarB5}, \ref{PDVVarB6}, \ref{PDVVarAB1}, \ref{PDVVarAB2}, and \ref{PDVVarAB3}. Therefore, following Definition~\ref{DnCat}, the $S$-sorted structure 
\[
\left(
\llbracket \mathbf{Pth}_{\boldsymbol{\mathcal{A}}^{(2)},s}\rrbracket,
\xi_{0,s}, \xi_{1,s}
\right)_{s\in S}
\]   
is an $S$-sorted $2$-category.
\end{proof}

The most important feature, as we will see immediately, is that $\llbracket \mathsf{Pth}_{\boldsymbol{\mathcal{A}}^{(2)}} \rrbracket$ is an $S$-sorted $2$-categorial $\Sigma$-algebra. To state this fact properly, we present the following propositions.

\begin{proposition}\label{PDVVarA7} Let $(\mathbf{s},s)$ be an element in $S^{\star}\times S$, $\sigma$ and operation symbol in $\Sigma_{\mathbf{s},s}$ and $(\llbracket\mathfrak{P}^{(2)}_{j} \rrbracket_{s_{j}})_{j\in\bb{\mathbf{s}}}$ a family of second-order path $\llbracket \cdot \rrbracket$-classes in $\llbracket \mathbf{Pth}_{\boldsymbol{\mathcal{A}}^{(2)}} \rrbracket_{\mathbf{s}}$, then the following equalities holds
\allowdisplaybreaks
\begin{multline*}
\sigma^{\llbracket \mathbf{Pth}_{\boldsymbol{\mathcal{A}}^{(2)}} \rrbracket}\left(
\left(
\mathrm{sc}^{0\llbracket \mathbf{Pth}_{\boldsymbol{\mathcal{A}}^{(2)}} \rrbracket}_{s_{j}}\left(
\Bigl\llbracket
\mathfrak{P}^{(2)}_{j}
\Bigr\rrbracket_{s_{j}}
\right)
\right)_{j\in\bb{\mathbf{s}}}
\right)
\\=
\mathrm{sc}^{0\llbracket \mathbf{Pth}_{\boldsymbol{\mathcal{A}}^{(2)}} \rrbracket}_{s}\left(
\sigma^{\llbracket \mathbf{Pth}_{\boldsymbol{\mathcal{A}}^{(2)}} \rrbracket}\left(
\left(
\Bigl\llbracket
\mathfrak{P}^{(2)}_{j}
\Bigr\rrbracket_{s_{j}}
\right)_{j\in\bb{\mathbf{s}}}
\right)
\right);
\end{multline*}
\allowdisplaybreaks
\begin{multline*}
\sigma^{\llbracket \mathbf{Pth}_{\boldsymbol{\mathcal{A}}^{(2)}} \rrbracket}\left(
\left(
\mathrm{tg}^{0\llbracket \mathbf{Pth}_{\boldsymbol{\mathcal{A}}^{(2)}} \rrbracket}_{s_{j}}\left(
\Bigl\llbracket
\mathfrak{P}^{(2)}_{j}
\Bigr\rrbracket_{s_{j}}
\right)
\right)_{j\in\bb{\mathbf{s}}}
\right)
\\=
\mathrm{tg}^{0\llbracket \mathbf{Pth}_{\boldsymbol{\mathcal{A}}^{(2)}} \rrbracket}_{s}\left(
\sigma^{\llbracket \mathbf{Pth}_{\boldsymbol{\mathcal{A}}^{(2)}} \rrbracket}\left(
\left(
\Bigl\llbracket
\mathfrak{P}^{(2)}_{j}
\Bigr\rrbracket_{s_{j}}
\right)_{j\in\bb{\mathbf{s}}}
\right)
\right).
\end{multline*}
\end{proposition}
\begin{proof}
We will only present the proof for the first equality. The other one is handled in a similar manner.

Note that the following chain of equalities holds
\begin{flushleft}
$\sigma^{\mathbf{Pth}_{\boldsymbol{\mathcal{A}}^{(2)}}}\left(
\left(
\mathrm{sc}^{0\mathbf{Pth}_{\boldsymbol{\mathcal{A}}^{(2)}}}_{s_{j}}\left(
\mathfrak{P}^{(2)}_{j}
\right)
\right)_{j\in\bb{\mathbf{s}}}
\right)$
\allowdisplaybreaks
\begin{align*}
&=
\sigma^{\mathbf{Pth}_{\boldsymbol{\mathcal{A}}^{(2)}}}\left(
\left(
\mathrm{ip}^{(2,0)\sharp}_{s_{j}}\left(
\mathrm{sc}^{(0,2)}_{s_{j}}\left(
\mathfrak{P}^{(2)}_{j}
\right)\right)
\right)_{j\in\bb{\mathbf{s}}}
\right)
\tag{1}
\\&=
\mathrm{ip}^{(2,0)\sharp}_{s}\left(
\mathrm{sc}^{(0,2)}_{s}\left(
\sigma^{\mathbf{Pth}_{\boldsymbol{\mathcal{A}}^{(2)}}}\left(
\left(
\mathfrak{P}^{(2)}_{j}
\right)_{j\in\bb{\mathbf{s}}}
\right)
\right)\right)
\tag{2}
\\&=
\mathrm{sc}^{0\mathbf{Pth}_{\boldsymbol{\mathcal{A}}^{(2)}}}_{s}\left(
\sigma^{\mathbf{Pth}_{\boldsymbol{\mathcal{A}}^{(2)}}}\left(
\left(
\mathfrak{P}^{(2)}_{j}
\right)_{j\in\bb{\mathbf{s}}}
\right)
\right).
\tag{3}
\end{align*}
\end{flushleft}

In the just stated chain of equalities, the first equality unravels the description of the $0$-source operation symbol in the many-sorted partial $\Sigma^{\boldsymbol{\mathcal{A}}^{(2)}}$-algebra $\mathbf{Pth}_{\boldsymbol{\mathcal{A}}^{(2)}}$ according to Proposition~\ref{PDPthCatAlg}; the second equality follows from the fact that $\mathrm{sc}^{(0,2)}$ and $\mathrm{ip}^{(2,0)\sharp}$ are $\Sigma$-homomorphisms according to Propositions~\ref{PDZHom} and~\ref{PDZHomIp}, respectively; finally, the last equality recovers the description of the $0$-source operation symbol in the many-sorted partial $\Sigma^{\boldsymbol{\mathcal{A}}^{(2)}}$-algebra $\mathbf{Pth}_{\boldsymbol{\mathcal{A}}^{(2)}}$ according to Proposition~\ref{PDPthCatAlg}.

Note that the desired equality follows from the fact that the $\llbracket\cdot\rrbracket$-classes of the two second-order paths above are equal.

This completes the proof.
\end{proof}

\begin{proposition}\label{PDVVarA8} Let $(\mathbf{s},s)$ be an element in $S^{\star}\times S$, $\sigma$ and operation symbol in $\Sigma_{\mathbf{s},s}$ and $(\llbracket\mathfrak{Q}^{(2)}_{j} \rrbracket_{s_{j}})_{j\in\bb{\mathbf{s}}}$ and $(\llbracket\mathfrak{P}^{(2)}_{j} \rrbracket_{s_{j}})_{j\in\bb{\mathbf{s}}}$ be two families of second-order path $\llbracket \cdot \rrbracket$-classes in $\llbracket \mathbf{Pth}_{\boldsymbol{\mathcal{A}}^{(2)}} \rrbracket_{\mathbf{s}}$ satisfying that, for every $j\in\bb{\mathbf{s}}$, 
\[
\mathrm{sc}^{0\llbracket \mathbf{Pth}_{\boldsymbol{\mathcal{A}}^{(2)}} \rrbracket}_{s_{j}}\left(
\Bigl\llbracket
\mathfrak{Q}^{(2)}
\Bigr\rrbracket_{s_{j}}
\right)
=
\mathrm{tg}^{0\llbracket \mathbf{Pth}_{\boldsymbol{\mathcal{A}}^{(2)}} \rrbracket}_{s_{j}}\left(
\Bigl\llbracket
\mathfrak{P}^{(2)}
\Bigr\rrbracket_{s_{j}}
\right),
\]
then the following equality holds
\allowdisplaybreaks
\begin{multline*}
\sigma^{\llbracket \mathbf{Pth}_{\boldsymbol{\mathcal{A}}^{(2)}} \rrbracket}\left(
\left(
\Bigl\llbracket
\mathfrak{Q}^{(2)}_{j}
\Bigr\rrbracket_{s_{j}}
\circ^{0\llbracket \mathbf{Pth}_{\boldsymbol{\mathcal{A}}^{(2)}} \rrbracket}_{s_{j}}
\Bigl\llbracket
\mathfrak{P}^{(2)}_{j}
\Bigr\rrbracket_{s_{j}}
\right)_{j\in\bb{\mathbf{s}}}
\right)
\\=
\sigma^{\llbracket \mathbf{Pth}_{\boldsymbol{\mathcal{A}}^{(2)}} \rrbracket}\left(
\left(
\Bigl\llbracket
\mathfrak{Q}^{(2)}_{j}
\Bigr\rrbracket_{s_{j}}
\right)_{j\in\bb{\mathbf{s}}}
\right)
\circ^{0\llbracket \mathbf{Pth}_{\boldsymbol{\mathcal{A}}^{(2)}} \rrbracket}_{s}
\sigma^{\llbracket \mathbf{Pth}_{\boldsymbol{\mathcal{A}}^{(2)}} \rrbracket}\left(
\left(
\Bigl\llbracket
\mathfrak{P}^{(2)}_{j}
\Bigr\rrbracket_{s_{j}}
\right)_{j\in\bb{\mathbf{s}}}
\right).
\end{multline*}
\end{proposition}
\begin{proof}
Let us recall from Definitions~\ref{DDUps} and~\ref{DDUpsCong} that the following equality holds 
\allowdisplaybreaks
\begin{multline*}
\left[
\sigma^{\mathbf{Pth}_{\boldsymbol{\mathcal{A}}^{(2)}}}\left(
\left(
\mathfrak{Q}^{(2)}_{j}
\circ^{0\mathbf{Pth}_{\boldsymbol{\mathcal{A}}^{(2)}}}_{s_{j}}
\mathfrak{P}^{(2)}_{j}
\right)_{j\in\bb{\mathbf{s}}}
\right)
\right]_{\Upsilon^{[1]}_{s}}
\\=
\left[
\sigma^{\mathbf{Pth}_{\boldsymbol{\mathcal{A}}^{(2)}}}\left(
\left(
\mathfrak{Q}^{(2)}_{j}
\right)_{j\in\bb{\mathbf{s}}}
\right)
\circ^{0\mathbf{Pth}_{\boldsymbol{\mathcal{A}}^{(2)}}}_{s}
\sigma^{\mathbf{Pth}_{\boldsymbol{\mathcal{A}}^{(2)}}}\left(
\left(
\mathfrak{P}^{(2)}_{j}
\right)_{j\in\bb{\mathbf{s}}}
\right)
\right]_{\Upsilon^{[1]}_{s}}
.
\end{multline*}

Therefore, taking $\llbracket \cdot \rrbracket$-classes, we conclude that 
\allowdisplaybreaks
\begin{multline*}
\left\llbracket
\sigma^{\mathbf{Pth}_{\boldsymbol{\mathcal{A}}^{(2)}}}\left(
\left(
\mathfrak{Q}^{(2)}_{j}
\circ^{0\mathbf{Pth}_{\boldsymbol{\mathcal{A}}^{(2)}}}_{s_{j}}
\mathfrak{P}^{(2)}_{j}
\right)_{j\in\bb{\mathbf{s}}}
\right)
\right\rrbracket_{s}
\\=
\left\llbracket
\sigma^{\mathbf{Pth}_{\boldsymbol{\mathcal{A}}^{(2)}}}\left(
\left(
\mathfrak{Q}^{(2)}_{j}
\right)_{j\in\bb{\mathbf{s}}}
\right)
\circ^{0\mathbf{Pth}_{\boldsymbol{\mathcal{A}}^{(2)}}}_{s}
\sigma^{\mathbf{Pth}_{\boldsymbol{\mathcal{A}}^{(2)}}}\left(
\left(
\mathfrak{P}^{(2)}_{j}
\right)_{j\in\bb{\mathbf{s}}}
\right)
\right\rrbracket_{s}
.
\end{multline*}

Taking into account the definition of the operation symbols in $\llbracket \mathbf{Pth}_{\boldsymbol{\mathcal{A}}^{(2)}} \rrbracket$, according to Proposition~\ref{PDVDCatAlg}, we conclude that 
\allowdisplaybreaks
\begin{multline*}
\sigma^{\llbracket \mathbf{Pth}_{\boldsymbol{\mathcal{A}}^{(2)}} \rrbracket}\left(
\left(
\Bigl\llbracket
\mathfrak{Q}^{(2)}_{j}
\Bigr\rrbracket_{s_{j}}
\circ^{0\llbracket \mathbf{Pth}_{\boldsymbol{\mathcal{A}}^{(2)}} \rrbracket}_{s_{j}}
\Bigl\llbracket
\mathfrak{P}^{(2)}_{j}
\Bigr\rrbracket_{s_{j}}
\right)_{j\in\bb{\mathbf{s}}}
\right)
\\=
\sigma^{\llbracket \mathbf{Pth}_{\boldsymbol{\mathcal{A}}^{(2)}} \rrbracket}\left(
\left(
\Bigl\llbracket
\mathfrak{Q}^{(2)}_{j}
\Bigr\rrbracket_{s_{j}}
\right)_{j\in\bb{\mathbf{s}}}
\right)
\circ^{0\llbracket \mathbf{Pth}_{\boldsymbol{\mathcal{A}}^{(2)}} \rrbracket}_{s}
\sigma^{\llbracket \mathbf{Pth}_{\boldsymbol{\mathcal{A}}^{(2)}} \rrbracket}\left(
\left(
\Bigl\llbracket
\mathfrak{P}^{(2)}_{j}
\Bigr\rrbracket_{s_{j}}
\right)_{j\in\bb{\mathbf{s}}}
\right).
\end{multline*}

This completes the proof.
\end{proof}

\begin{proposition}\label{PDVVarB7} Let $(\mathbf{s},s)$ be an element in $S^{\star}\times S$, $\sigma$ and operation symbol in $\Sigma_{\mathbf{s},s}$ and $(\llbracket\mathfrak{P}^{(2)}_{j} \rrbracket_{s_{j}})_{j\in\bb{\mathbf{s}}}$ a family of second-order path $\llbracket \cdot \rrbracket$-classes in $\llbracket \mathbf{Pth}_{\boldsymbol{\mathcal{A}}^{(2)}} \rrbracket_{\mathbf{s}}$, then the following equalities holds
\allowdisplaybreaks
\begin{multline*}
\sigma^{\llbracket \mathbf{Pth}_{\boldsymbol{\mathcal{A}}^{(2)}} \rrbracket}\left(
\left(
\mathrm{sc}^{1\llbracket \mathbf{Pth}_{\boldsymbol{\mathcal{A}}^{(2)}} \rrbracket}_{s_{j}}\left(
\Bigl\llbracket
\mathfrak{P}^{(2)}_{j}
\Bigr\rrbracket_{s_{j}}
\right)
\right)_{j\in\bb{\mathbf{s}}}
\right)
\\=
\mathrm{sc}^{1\llbracket \mathbf{Pth}_{\boldsymbol{\mathcal{A}}^{(2)}} \rrbracket}_{s}\left(
\sigma^{\llbracket \mathbf{Pth}_{\boldsymbol{\mathcal{A}}^{(2)}} \rrbracket}\left(
\left(
\Bigl\llbracket
\mathfrak{P}^{(2)}_{j}
\Bigr\rrbracket_{s_{j}}
\right)_{j\in\bb{\mathbf{s}}}
\right)
\right);
\end{multline*}
\allowdisplaybreaks
\begin{multline*}
\sigma^{\llbracket \mathbf{Pth}_{\boldsymbol{\mathcal{A}}^{(2)}} \rrbracket}\left(
\left(
\mathrm{tg}^{1\llbracket \mathbf{Pth}_{\boldsymbol{\mathcal{A}}^{(2)}} \rrbracket}_{s_{j}}\left(
\Bigl\llbracket
\mathfrak{P}^{(2)}_{j}
\Bigr\rrbracket_{s_{j}}
\right)
\right)_{j\in\bb{\mathbf{s}}}
\right)
\\=
\mathrm{tg}^{1\llbracket \mathbf{Pth}_{\boldsymbol{\mathcal{A}}^{(2)}} \rrbracket}_{s}\left(
\sigma^{\llbracket \mathbf{Pth}_{\boldsymbol{\mathcal{A}}^{(2)}} \rrbracket}\left(
\left(
\Bigl\llbracket
\mathfrak{P}^{(2)}_{j}
\Bigr\rrbracket_{s_{j}}
\right)_{j\in\bb{\mathbf{s}}}
\right)
\right).
\end{multline*}
\end{proposition}
\begin{proof}
We will only present the proof for the first equality. The other one is handled in a similar manner.

Note that the following chain of equalities holds
\begin{flushleft}
$\sigma^{\mathbf{Pth}_{\boldsymbol{\mathcal{A}}^{(2)}}}\left(
\left(
\mathrm{sc}^{1\mathbf{Pth}_{\boldsymbol{\mathcal{A}}^{(2)}}}_{s_{j}}\left(
\mathfrak{P}^{(2)}_{j}
\right)
\right)_{j\in\bb{\mathbf{s}}}
\right)$
\allowdisplaybreaks
\begin{align*}
\qquad&=
\sigma^{\mathbf{Pth}_{\boldsymbol{\mathcal{A}}^{(2)}}}\left(
\left(
\mathrm{ip}^{(2,[1])\sharp}_{s_{j}}\left(
\mathrm{sc}^{([1],2)}_{s_{j}}\left(
\mathfrak{P}^{(2)}_{j}
\right)\right)
\right)_{j\in\bb{\mathbf{s}}}
\right)
\tag{1}
\\&=
\mathrm{ip}^{(2,[1])\sharp}_{s}\left(
\mathrm{sc}^{([1],2)}_{s}\left(
\sigma^{\mathbf{Pth}_{\boldsymbol{\mathcal{A}}^{(2)}}}\left(
\left(
\mathfrak{P}^{(2)}_{j}
\right)_{j\in\bb{\mathbf{s}}}
\right)
\right)\right)
\tag{2}
\\&=
\mathrm{sc}^{1\mathbf{Pth}_{\boldsymbol{\mathcal{A}}^{(2)}}}_{s}\left(
\sigma^{\mathbf{Pth}_{\boldsymbol{\mathcal{A}}^{(2)}}}\left(
\left(
\mathfrak{P}^{(2)}_{j}
\right)_{j\in\bb{\mathbf{s}}}
\right)
\right).
\tag{3}
\end{align*}
\end{flushleft}

In the just stated chain of equalities, the first equality unravels the description of the $1$-source operation symbol in the many-sorted partial $\Sigma^{\boldsymbol{\mathcal{A}}^{(2)}}$-algebra $\mathbf{Pth}_{\boldsymbol{\mathcal{A}}^{(2)}}$ according to Proposition~\ref{PDPthDCatAlg}; the second equality follows from the fact that $\mathrm{sc}^{([1],2)}$ and $\mathrm{ip}^{(2,[1])\sharp}$ are $\Sigma^{\boldsymbol{\mathcal{A}}}$-homomorphisms according to Propositions~\ref{PDUCatHom} and~\ref{PDUIpCatHom}, respectively; finally, the last equality recovers the description of the $1$-source operation symbol in the many-sorted partial $\Sigma^{\boldsymbol{\mathcal{A}}^{(2)}}$-algebra $\mathbf{Pth}_{\boldsymbol{\mathcal{A}}^{(2)}}$ according to Proposition~\ref{PDPthDCatAlg}.

Note that the desired equality follows from the fact that the $\llbracket\cdot\rrbracket$-classes of the two second-order paths above are equal.

This completes the proof.
\end{proof}

\begin{proposition}\label{PDVVarB8} Let $(\mathbf{s},s)$ be an element in $S^{\star}\times S$, $\sigma$ and operation symbol in $\Sigma_{\mathbf{s},s}$ and $(\llbracket\mathfrak{Q}^{(2)}_{j} \rrbracket_{s_{j}})_{j\in\bb{\mathbf{s}}}$ and $(\llbracket\mathfrak{P}^{(2)}_{j} \rrbracket_{s_{j}})_{j\in\bb{\mathbf{s}}}$ be two families of second-order path $\llbracket \cdot \rrbracket$-classes in $\llbracket \mathbf{Pth}_{\boldsymbol{\mathcal{A}}^{(2)}} \rrbracket_{\mathbf{s}}$ satisfying that, for every $j\in\bb{\mathbf{s}}$, 
\[
\mathrm{sc}^{1\llbracket \mathbf{Pth}_{\boldsymbol{\mathcal{A}}^{(2)}} \rrbracket}_{s_{j}}\left(
\Bigl\llbracket
\mathfrak{Q}^{(2)}
\Bigr\rrbracket_{s_{j}}
\right)
=
\mathrm{tg}^{1\llbracket \mathbf{Pth}_{\boldsymbol{\mathcal{A}}^{(2)}} \rrbracket}_{s_{j}}\left(
\Bigl\llbracket
\mathfrak{P}^{(2)}
\Bigr\rrbracket_{s_{j}}
\right),
\]
then the following equality holds
\allowdisplaybreaks
\begin{multline*}
\sigma^{\llbracket \mathbf{Pth}_{\boldsymbol{\mathcal{A}}^{(2)}} \rrbracket}\left(
\left(
\Bigl\llbracket
\mathfrak{Q}^{(2)}_{j}
\Bigr\rrbracket_{s_{j}}
\circ^{1\llbracket \mathbf{Pth}_{\boldsymbol{\mathcal{A}}^{(2)}} \rrbracket}_{s_{j}}
\Bigl\llbracket
\mathfrak{P}^{(2)}_{j}
\Bigr\rrbracket_{s_{j}}
\right)_{j\in\bb{\mathbf{s}}}
\right)
\\=
\sigma^{\llbracket \mathbf{Pth}_{\boldsymbol{\mathcal{A}}^{(2)}} \rrbracket}\left(
\left(
\Bigl\llbracket
\mathfrak{Q}^{(2)}_{j}
\Bigr\rrbracket_{s_{j}}
\right)_{j\in\bb{\mathbf{s}}}
\right)
\circ^{1\llbracket \mathbf{Pth}_{\boldsymbol{\mathcal{A}}^{(2)}} \rrbracket}_{s}
\sigma^{\llbracket \mathbf{Pth}_{\boldsymbol{\mathcal{A}}^{(2)}} \rrbracket}\left(
\left(
\Bigl\llbracket
\mathfrak{P}^{(2)}_{j}
\Bigr\rrbracket_{s_{j}}
\right)_{j\in\bb{\mathbf{s}}}
\right).
\end{multline*}
\end{proposition}
\begin{proof} 
We will consider two cases according to the nature of the families of second-order paths $(\mathfrak{Q}^{(2)}_{j})_{j\in\bb{\mathbf{s}}}$ and $(\mathfrak{P}^{(2)}_{j})_{j\in\bb{\mathbf{s}}}$. It could be the case that (1) all the second-order paths above are $(2,[1])$-identity second-order paths or (2) there exists at least one of them that is not a $(2,[1])$-identity second-order path.

For the first case, assume that, for suitable families of path terms $([Q_{j}]_{s_{j}})_{j\in\bb{\mathbf{s}}}$ and $([P_{j}]_{s_{j}})_{j\in\bb{\mathbf{s}}}$ in $[\mathrm{PT}_{\boldsymbol{\mathcal{A}}}]_{\mathbf{s}}$, it is the case that, for every $j\in\bb{\mathbf{s}}$, 
\begin{align*}
\mathfrak{Q}^{(2)}_{j}&=
\mathrm{ip}^{(2,[1])\sharp}_{s_{j}}\left(
[Q_{j}]_{s_{j}}
\right);
&
\mathfrak{P}^{(2)}_{j}&=
\mathrm{ip}^{(2,[1])\sharp}_{s_{j}}\left(
[P_{j}]_{s_{j}}
\right).
\end{align*}

Note that, in order to satisfy the equality that holds by hypothesis, for every $j\in\bb{\mathbf{s}}$, the following chain of equalities must hold
\begin{align*}
\Bigl\llbracket \mathrm{ip}^{(2,[1])\sharp}_{s_{j}}\left(
\mathrm{sc}^{([1],2)}_{s_{j}}\left(
\mathfrak{Q}^{(2)}_{j}
\right)
\right) \Bigr\rrbracket_{s_{j}}
&=
\Bigl\llbracket 
\mathrm{sc}^{1\mathbf{Pth}_{\boldsymbol{\mathcal{A}}^{(2)}}}_{s_{j}}\left(
\mathfrak{Q}^{(2)}_{j}
\right)
\Bigr\rrbracket_{s_{j}}
\tag{1}
\\
&=
\mathrm{sc}^{1\llbracket \mathbf{Pth}_{\boldsymbol{\mathcal{A}}^{(2)}} \rrbracket}_{s_{j}}\left(
\llbracket
\mathfrak{Q}^{(2)}_{j}
\rrbracket
\right)
\tag{2}
\\&=
\mathrm{tg}^{1\llbracket \mathbf{Pth}_{\boldsymbol{\mathcal{A}}^{(2)}} \rrbracket}_{s_{j}}\left(
\llbracket
\mathfrak{P}^{(2)}_{j}
\rrbracket
\right)
\tag{3}
\\&=
\Bigl\llbracket 
\mathrm{tg}^{1\mathbf{Pth}_{\boldsymbol{\mathcal{A}}^{(2)}}}_{s_{j}}\left(
\mathfrak{P}^{(2)}_{j}
\right)
\Bigr\rrbracket_{s_{j}}
\tag{4}
\\&=
\Bigl\llbracket \mathrm{ip}^{(2,[1])\sharp}_{s_{j}}\left(
\mathrm{tg}^{([1],2)}_{s_{j}}\left(
\mathfrak{P}^{(2)}
\right)
\right) \Bigr\rrbracket_{s_{j}}.
\tag{5}
\end{align*}

In the just stated chain of equalities, the first equality recovers the description of the $1$-source operation symbol in the many-sorted partial $\Sigma^{\boldsymbol{\mathcal{A}}^{(2)}}$-algebra $\mathbf{Pth}_{\boldsymbol{\mathcal{A}}^{(2)}}$, according to Proposition~\ref{PDPthDCatAlg}; the second equality recovers the description of the $1$-source operation symbol in the many-sorted partial $\Sigma^{\boldsymbol{\mathcal{A}}^{(2)}}$-algebra $\llbracket \mathbf{Pth}_{\boldsymbol{\mathcal{A}}^{(2)}}\rrbracket $, according to Proposition~\ref{PDVDCatAlg}; the third equality holds by assumption; the fourth equality unravels the description of the $1$-target operation symbol in the many-sorted partial $\Sigma^{\boldsymbol{\mathcal{A}}^{(2)}}$-algebra $\llbracket \mathbf{Pth}_{\boldsymbol{\mathcal{A}}^{(2)}}\rrbracket $, according to Proposition~\ref{PDVDCatAlg}; finally, the last equality recovers the description of the $1$-target operation symbol in the many-sorted partial $\Sigma^{\boldsymbol{\mathcal{A}}^{(2)}}$-algebra $\mathbf{Pth}_{\boldsymbol{\mathcal{A}}^{(2)}}$, according to Proposition~\ref{PDPthDCatAlg}.

Taking into account the above equality and Proposition~\ref{PDVDUId}, we conclude that, for every $j\in\bb{\mathbf{s}}$, 
\[
\mathrm{ip}^{(2,[1])\sharp}_{s_{j}}\left(
\mathrm{sc}^{([1],2)}_{s_{j}}\left(
\mathfrak{Q}^{(2)}_{j}
\right)
\right)
=
\mathrm{ip}^{(2,[1])\sharp}_{s_{j}}\left(
\mathrm{tg}^{([1],2)}_{s_{j}}\left(
\mathfrak{P}^{(2)}_{j}
\right)
\right).
\]

Applying the $([1],2)$-source mapping on both sides above we conclude, according to Proposition~\ref{PDBasicEq} that, for every $j\in\bb{\mathbf{s}}$,
\[
[Q_{j}]_{s_{j}}=
\mathrm{sc}^{([1],2)}_{s_{j}}\left(
\mathfrak{Q}^{(2)}_{j}
\right)
=
\mathrm{tg}^{([1],2)}_{s_{j}}\left(
\mathfrak{P}^{(2)}_{j}
\right)
=
[P_{j}]_{s_{j}}.
\]

Consequently, for every $j\in\bb{\mathbf{s}}$, the two second-order paths $\mathfrak{Q}^{(2)}_{j}$ and $\mathfrak{P}^{(2)}_{j}$ are equal, since
\[
\mathfrak{Q}^{(2)}_{j}
=\mathrm{ip}^{(2,[1])\sharp}_{s_{j}}\left([Q_{j}]_{s_{j}}\right)
=\mathrm{ip}^{(2,[1])\sharp}_{s_{j}}\left([P_{j}]_{s_{j}}\right)
=\mathfrak{P}^{(2)}_{j}.
\]

Then the following chain of equalities holds
\begin{flushleft}
$\sigma^{\mathbf{Pth}_{\boldsymbol{\mathcal{A}}^{(2)}}}
\left(
\left(
\mathfrak{Q}^{(2)}_{j}
\circ^{1\mathbf{Pth}_{\boldsymbol{\mathcal{A}}^{(2)}}}_{s_{j}}
\mathfrak{P}^{(2)}_{j}
\right)_{j\in\bb{\mathbf{s}}}
\right)
$
\allowdisplaybreaks
\begin{align*}
&=
\sigma^{\mathbf{Pth}_{\boldsymbol{\mathcal{A}}^{(2)}}}
\left(
\left(
\mathrm{ip}^{(2,[1])\sharp}_{s_{j}}\left(
[P]_{s_{j}}
\right)
\circ^{1\mathbf{Pth}_{\boldsymbol{\mathcal{A}}^{(2)}}}_{s_{j}}
\mathrm{ip}^{(2,[1])\sharp}_{s_{j}}\left(
[P]_{s_{j}}
\right)
\right)_{j\in\bb{\mathbf{s}}}
\right)
\tag{1}
\\&=
\sigma^{\mathbf{Pth}_{\boldsymbol{\mathcal{A}}^{(2)}}}
\left(
\left(
\mathrm{ip}^{(2,[1])\sharp}_{s_{j}}\left(
[P]_{s_{j}}
\right)
\right)_{j\in\bb{\mathbf{s}}}
\right)
\tag{2}
\\&=
\mathrm{ip}^{(2,[1])\sharp}_{s}
\left(
\sigma^{[\mathbf{PT}_{\boldsymbol{\mathcal{A}}}]}
\left(
\left(
[P]_{s_{j}}
\right)_{j\in\bb{\mathbf{s}}}
\right)
\right)
\tag{3}
\\&=
\mathrm{ip}^{(2,[1])\sharp}_{s}
\left(
\sigma^{[\mathbf{PT}_{\boldsymbol{\mathcal{A}}}]}
\left(
\left(
[P]_{s_{j}}
\right)_{j\in\bb{\mathbf{s}}}
\right)
\right)
\circ^{1\mathbf{Pth}_{\boldsymbol{\mathcal{A}}^{(2)}}}_{s}
\mathrm{ip}^{(2,[1])\sharp}_{s}
\left(
\sigma^{[\mathbf{PT}_{\boldsymbol{\mathcal{A}}}]}
\left(
\left(
[P]_{s_{j}}
\right)_{j\in\bb{\mathbf{s}}}
\right)
\right)
\tag{4}
\\&=
\sigma^{\mathbf{Pth}_{\boldsymbol{\mathcal{A}}^{(2)}}}\left(
\left(
\mathrm{ip}^{(2,[1])\sharp}_{s_{j}}\left(
[P]_{s_{j}}
\right)
\right)_{j\in\bb{\mathbf{s}}}
\right)
\circ^{1\mathbf{Pth}_{\boldsymbol{\mathcal{A}}^{(2)}}}_{s}
\\&\qquad\qquad\qquad\qquad\qquad\qquad\qquad\qquad\qquad
\sigma^{\mathbf{Pth}_{\boldsymbol{\mathcal{A}}^{(2)}}}\left(
\left(
\mathrm{ip}^{(2,[1])\sharp}_{s_{j}}\left(
[P]_{s_{j}}
\right)
\right)_{j\in\bb{\mathbf{s}}}
\right)
\tag{5}
\\&=
\sigma^{\mathbf{Pth}_{\boldsymbol{\mathcal{A}}^{(2)}}}\left(
\left(
\mathfrak{Q}^{(2)}_{j}
\right)_{j\in\bb{\mathbf{s}}}
\right)
\circ^{1\mathbf{Pth}_{\boldsymbol{\mathcal{A}}^{(2)}}}_{s}
\sigma^{\mathbf{Pth}_{\boldsymbol{\mathcal{A}}^{(2)}}}\left(
\left(
\mathfrak{P}^{(2)}_{j}
\right)_{j\in\bb{\mathbf{s}}}
\right).
\tag{6}
\end{align*}
\end{flushleft}

In the just stated chain of equalities, the first equality follows from the description of the second-order paths as $(2,[1])$-identity second-order paths; the second equality follows from the fact that, according to Definition~\ref{DDPthComp}, the $(2,[1])$-identity second-order paths are idempotent for the $1$-composition; the third equality follows from Proposition~\ref{PDUSigma};  the fourth equality follows from the fact that, according to Definition~\ref{DDPthComp}, the $(2,[1])$-identity second-order paths are idempotent for the $1$-composition; the fifth equality follows from Proposition~\ref{PDUSigma}; finally, the last equality recovers the  description of the second-order paths as $(2,[1])$-identity second-order paths.

Note that the desired equality follows from the fact that the $\llbracket\cdot\rrbracket$-classes of the two second-order paths above are equal.

This completes case~(1), where all second-order paths above are $(2,[1])$-identity second-order paths.

Now, consider the case~(2), i.e., the case in which at least one of the second-order paths under consideration is not a $(2,[1])$-identity second-order path.

Consider, on one hand, the second-order path 
\[
\sigma^{\mathbf{Pth}_{\boldsymbol{\mathcal{A}}^{(2)}} }\left(
\left(
\mathfrak{Q}^{(2)}_{j}
\circ^{1\mathbf{Pth}_{\boldsymbol{\mathcal{A}}^{(2)}} }_{s_{j}}
\mathfrak{P}^{(2)}_{j}
\right)_{j\in\bb{\mathbf{s}}}
\right).
\]

An schematic representation of this second-order path is included in Figure~\ref{FDVVarB8Lft}. Note that we have removed the superscripts $[\mathbf{PT}_{\boldsymbol{\mathcal{A}}}]$ for clarity.  This is, in virtue of Corollary~\ref{CDPthWB}, a coherent head-constant echelonless second-order path associated to the $\sigma$ operation symbol. Note that the second-order path extraction procedure from Lemma~\ref{LDPthExtract} applied to it, retrieves the family of second-order paths
\[
\left(
\mathfrak{Q}^{(2)}_{j}
\circ^{1\mathbf{Pth}_{\boldsymbol{\mathcal{A}}^{(2)}}}_{s_{j}}
\mathfrak{P}^{(2)}_{j}
\right)_{j\in\bb{\mathbf{s}}}.
\] 

Therefore, the second-order Curry-Howard mapping at this second-order path is given by the following expression.
\allowdisplaybreaks
\begin{multline*}
\mathrm{CH}^{(2)}_{s}\left(
\sigma^{\mathbf{Pth}_{\boldsymbol{\mathcal{A}}^{(2)}} }\left(
\left(
\mathfrak{Q}^{(2)}_{j}
\circ^{1\mathbf{Pth}_{\boldsymbol{\mathcal{A}}^{(2)}} }_{s_{j}}
\mathfrak{P}^{(2)}_{j}
\right)_{j\in\bb{\mathbf{s}}}
\right)
\right)
\\=
\sigma^{\mathbf{T}_{\Sigma^{\boldsymbol{\mathcal{A}}^{(2)}}}(X)}\left(
\left(
\mathrm{CH}^{(2)}_{s_{j}}\left(
\mathfrak{Q}^{(2)}_{j}
\circ^{1\mathbf{Pth}_{\boldsymbol{\mathcal{A}}^{(2)}} }_{s_{j}}
\mathfrak{P}^{(2)}_{j}
\right)
\right)_{j\in\bb{\mathbf{s}}}
\right).
\end{multline*}

\begin{figure}
\centering
\begin{tikzpicture}[
fletxa/.style={thick, ->}]
\tikzstyle{every state}=[very thick, draw=blue!50,fill=blue!20, inner sep=2pt,minimum size=20pt]

\node[] (a0) at (1.15,0) [] {$\sigma\Biggl($};
\node[] (a1) at (2.75,0) [] {$\mathrm{sc}^{([1],2)}_{s_{0}}\left(\mathfrak{P}^{(2)}_{0}\right)$}; 
\node[] (a2) at (4.5,0) [] {$,\cdots ,$}; 
\node[] (a3) at (6.4,0) [] {$\mathrm{sc}^{([1],2)}_{s_{\bb{\mathbf{s}}-1}}\left(\mathfrak{P}^{(2)}_{\bb{\mathbf{s}}-1}\right)$};
\node[] (a4) at (8,0) [] {$\Biggr)$};

\node[] (b1) at (3.05,-.95) [] {\rotatebox{-90}{
$\xymatrix@C=30pt{
\,
\ar@{=>}[r]^-{\text{
\rotatebox{90}{$\mathfrak{P}^{(2)}_{0}$}
}}
&
\,
}$}}; 
\node[] (b3) at (6.4,-1) [] {\rotatebox{-90}{$\,=$}};

\node[] (c0) at (1.15,-2) [] {$\sigma\Biggl($};
\node[] (c1) at (2.75,-2) [] {$\mathrm{tg}^{([1],2)}_{s_{0}}\left(\mathfrak{P}^{(2)}_{0}\right)$}; 
\node[] (c2) at (4.5,-2) [] {$,\cdots ,$}; 
\node[] (c3) at (6.4,-2) [] {$\vdots$};
\node[] (c4) at (8,-2) [] {$\Biggr)$};

\node[] (d1) at (2.75,-3) [] {\rotatebox{-90}{$\,=$}};  
\node[] (d3) at (6.4,-3) [] {\rotatebox{-90}{$\,=$}};

\node[] (e0) at (1.15,-4) [] {$\sigma\Biggl($};
\node[] (e1) at (2.75,-4) [] {$\mathrm{sc}^{([1],2)}_{s_{0}}\left(\mathfrak{Q}^{(2)}_{0}\right)$}; 
\node[] (e2) at (4.5,-4) [] {$,\cdots ,$}; 
\node[] (e3) at (6.4,-4) [] {$\vdots$};
\node[] (e4) at (8,-4) [] {$\Biggr)$};

\node[] (f1) at (3.05,-4.95) [] {\rotatebox{-90}{
$\xymatrix@C=30pt{
\,
\ar@{=>}[r]^-{\text{
\rotatebox{90}{$\mathfrak{Q}^{(2)}_{0}$}
}}
&
\,
}$}}; 
\node[] (f3) at (6.4,-5) [] {\rotatebox{-90}{$\,=$}};

\node[] (g0) at (1.15,-6) [] {$\sigma\Biggl($};
\node[] (g1) at (2.75,-6) [] {$\mathrm{tg}^{([1],2)}_{s_{0}}\left(\mathfrak{Q}^{(2)}_{0}\right)$}; 
\node[] (g2) at (4.5,-6) [] {$,\cdots ,$}; 
\node[] (g3) at (6.4,-6) [] {$\vdots$};
\node[] (g4) at (8,-6) [] {$\Biggr)$};

\node[] (h1) at (2.75,-7) [] {\color{white}{$\mathrm{tg}^{([1],2)}_{s_{0}}\left(\mathfrak{Q}^{(2)}_{0}\right)$}}; 
\node[] (h1b) at (2.75,-7) [] {\rotatebox{-90}{$\,=$}}; 
\node[] (h2) at (4.5,-7) [] {$\ddots$}; 
\node[] (h3) at (6.4,-7) [] {\color{white}{$\mathrm{sc}^{([1],2)}_{s_{\bb{\mathbf{s}}-1}}\left(\mathfrak{P}^{(2)}_{\bb{\mathbf{s}}-1}\right)$}};
\node[] (h3b) at (6.4,-7) [] {\rotatebox{-90}{$\,=$}}; 

\node[] (i0) at (1.15,-8) [] {$\sigma\Biggl($};
\node[] (i1) at (2.75,-8) [] {$\vdots$}; 
\node[] (i2) at (4.5,-8) [] {$,\cdots ,$}; 
\node[] (i3) at (6.4,-8) [] {$\mathrm{sc}^{([1],2)}_{s_{\bb{\mathbf{s}}-1}}\left(\mathfrak{P}^{(2)}_{\bb{\mathbf{s}}-1}\right)$};
\node[] (i4) at (8,-8) [] {$\Biggr)$};

\node[] (j1) at (2.75,-9) [] {\rotatebox{-90}{$\,=$}}; 
\node[] (j3) at (6.85,-8.95) []  {\rotatebox{-90}{
$\xymatrix@C=30pt{
\,
\ar@{=>}[r]^-{\text{
\rotatebox{90}{$\mathfrak{P}^{(2)}_{\bb{\mathbf{s}}-1}$}
}}
&
\,
}$}};

\node[] (k0) at (1.15,-10) [] {$\sigma\Biggl($};
\node[] (k1) at (2.75,-10) [] {$\vdots$}; 
\node[] (k2) at (4.5,-10) [] {$,\cdots ,$}; 
\node[] (k3) at (6.4,-10) [] {$\mathrm{tg}^{([1],2)}_{s_{\bb{\mathbf{s}}-1}}\left(\mathfrak{P}^{(2)}_{\bb{\mathbf{s}}-1}\right)$};
\node[] (k4) at (8,-10) [] {$\Biggr)$};

\node[] (l1) at (2.75,-11) [] {\rotatebox{-90}{$\,\cdots$}}; 
\node[] (l3) at (6.4,-11) []   {\rotatebox{-90}{$\,=$}}; 

\node[] (m0) at (1.15,-12) [] {$\sigma\Biggl($};
\node[] (m1) at (2.75,-12) [] {$\vdots$}; 
\node[] (m2) at (4.5,-12) [] {$,\cdots ,$}; 
\node[] (m3) at (6.4,-12) [] {$\mathrm{sc}^{([1],2)}_{s_{\bb{\mathbf{s}}-1}}\left(\mathfrak{Q}^{(2)}_{\bb{\mathbf{s}}-1}\right)$};
\node[] (m4) at (8,-12) [] {$\Biggr)$};

\node[] (n1) at (2.75,-13) [] {\rotatebox{-90}{$\,=$}}; 
\node[] (n3) at (6.85,-12.95) []  {\rotatebox{-90}{
$\xymatrix@C=30pt{
\,
\ar@{=>}[r]^-{\text{
\rotatebox{90}{$\mathfrak{Q}^{(2)}_{\bb{\mathbf{s}}-1}$}
}}
&
\,
}$}};

\node[] (o0) at (1.15,-14) [] {$\sigma\Biggl($};
\node[] (o1) at (2.75,-14) [] {$\mathrm{tg}^{([1],2)}_{s_{\bb{\mathbf{s}}-1}}\left(\mathfrak{Q}^{(2)}_{s_{0}}\right)$}; 
\node[] (o2) at (4.5,-14) [] {$,\cdots ,$}; 
\node[] (o3) at (6.4,-14) [] {$\mathrm{tg}^{([1],2)}_{s_{\bb{\mathbf{s}}-1}}\left(\mathfrak{Q}^{(2)}_{\bb{\mathbf{s}}-1}\right)$};
\node[] (o4) at (8,-14) [] {$\Biggr)$};

\tikzset{encercla/.style={draw=black, line width=.5pt, inner sep=0pt, rectangle, rounded corners}};
\node [encercla,  fit=(a1) ] {} ;
\node [encercla,  fit=(a3)(i3) ] {} ;
\node [encercla,  fit=(c1)(d1)(e1) ] {} ;
\node [encercla,  fit=(g1)(o1) ] {} ;
\node [encercla,  fit=(k3)(m3) ] {} ;
\node [encercla,  fit=(o3) ] {} ;
\tikzset{encercla2/.style={draw=white, dashed, line width=.65pt, inner sep=0pt, rectangle, rounded corners}};
\node [encercla2,  fit=(h1) ] {} ;
\node [encercla2,  fit=(h3) ] {} ;
\end{tikzpicture}
\caption{
$\sigma^{\mathbf{Pth}_{\boldsymbol{\mathcal{A}}^{(2)}} }(
(
\mathfrak{Q}^{(2)}_{j}
\circ^{1\mathbf{Pth}_{\boldsymbol{\mathcal{A}}^{(2)}} }_{s_{j}}
\mathfrak{P}^{(2)}_{j}
)_{j\in\bb{\mathbf{s}}}
)$.
}\label{FDVVarB8Lft}
\end{figure}

Now consider, on the other hand, the second-order path
\[
\sigma^{\mathbf{Pth}_{\boldsymbol{\mathcal{A}}^{(2)}}}\left(
\left(
\mathfrak{Q}^{(2)}_{j}
\right)_{j\in \bb{\mathbf{s}}}
\right)
\circ^{1\mathbf{Pth}_{\boldsymbol{\mathcal{A}}^{(2)}}}
\sigma^{\mathbf{Pth}_{\boldsymbol{\mathcal{A}}^{(2)}}}\left(
\left(
\mathfrak{P}^{(2)}_{j}
\right)_{j\in \bb{\mathbf{s}}}
\right).
\]

An schematic representation of this second-order path is included in Figure~\ref{FDVVarB8Rgt}. Note that we have removed the superscripts $[\mathbf{PT}_{\boldsymbol{\mathcal{A}}}]$ for clarity.  This is, in virtue of Corollary~\ref{CDPthWB}, a coherent head-constant echelonless second-order path associated to the $\sigma$ operation symbol. Note that the second-order path extraction procedure from Lemma~\ref{LDPthExtract} applied to it, retrieves the family of second-order paths
\[
\left(
\mathfrak{Q}^{(2)}_{j}
\circ^{1\mathbf{Pth}_{\boldsymbol{\mathcal{A}}^{(2)}}}_{s_{j}}
\mathfrak{P}^{(2)}_{j}
\right)_{j\in\bb{\mathbf{s}}}.
\] 

Therefore, the second-order Curry-Howard mapping at this second-order path is given by the following expression.
\allowdisplaybreaks
\begin{multline*}
\mathrm{CH}^{(2)}_{s}\left(
\sigma^{\mathbf{Pth}_{\boldsymbol{\mathcal{A}}^{(2)}}}\left(
\left(
\mathfrak{Q}^{(2)}_{j}
\right)_{j\in \bb{\mathbf{s}}}
\right)
\circ^{1\mathbf{Pth}_{\boldsymbol{\mathcal{A}}^{(2)}}}
\sigma^{\mathbf{Pth}_{\boldsymbol{\mathcal{A}}^{(2)}}}\left(
\left(
\mathfrak{P}^{(2)}_{j}
\right)_{j\in \bb{\mathbf{s}}}
\right)
\right)
\\=
\sigma^{\mathbf{T}_{\Sigma^{\boldsymbol{\mathcal{A}}^{(2)}}}(X)}\left(
\left(
\mathrm{CH}^{(2)}_{s_{j}}\left(
\mathfrak{Q}^{(2)}_{j}
\circ^{1\mathbf{Pth}_{\boldsymbol{\mathcal{A}}^{(2)}} }_{s_{j}}
\mathfrak{P}^{(2)}_{j}
\right)
\right)_{j\in\bb{\mathbf{s}}}
\right).
\end{multline*}

\begin{figure}
\centering
\begin{tikzpicture}[
fletxa/.style={thick, ->}]
\tikzstyle{every state}=[very thick, draw=blue!50,fill=blue!20, inner sep=2pt,minimum size=20pt]

\node[] (a0) at (1.15,0) [] {$\sigma \Biggl($};
\node[] (a1) at (2.75,0) [] {$\mathrm{sc}^{([1],2)}_{s_{0}}\left(\mathfrak{P}^{(2)}_{0}\right)$}; 
\node[] (a2) at (4.5,0) [] {$,\cdots ,$}; 
\node[] (a3) at (6.4,0) [] {$\mathrm{sc}^{([1],2)}_{s_{\bb{\mathbf{s}}-1}}\left(\mathfrak{P}^{(2)}_{\bb{\mathbf{s}}-1}\right)$};
\node[] (a4) at (8,0) [] {$\Biggr)$};

\node[] (b1) at (3.05,-.95) [] {\rotatebox{-90}{
$\xymatrix@C=30pt{
\,
\ar@{=>}[r]^-{\text{
\rotatebox{90}{$\mathfrak{P}^{(2)}_{0}$}
}}
&
\,
}$}}; 
\node[] (b3) at (6.4,-1) [] {\rotatebox{-90}{$\,=$}};

\node[] (c0) at (1.15,-2) [] {$\sigma \Biggl($};
\node[] (c1) at (2.75,-2) [] {$\mathrm{tg}^{([1],2)}_{s_{0}}\left(\mathfrak{P}^{(2)}_{0}\right)$}; 
\node[] (c2) at (4.5,-2) [] {$,\cdots ,$}; 
\node[] (c3) at (6.4,-2) [] {$\vdots$};
\node[] (c4) at (8,-2) [] {$\Biggr)$};

\node[] (d1) at (2.75,-3) [] {\color{white}{$\mathrm{tg}^{([1],2)}_{s_{0}}\left(\mathfrak{P}^{(2)}_{0}\right)$}}; 
\node[] (d1b) at (2.75,-3) [] {\rotatebox{-90}{$\,=$}}; 
\node[] (d2) at (4.5,-3) [] {$\ddots$}; 
\node[] (d3) at (6.4,-3) [] {\color{white}{$\mathrm{sc}^{([1],2)}_{s_{\bb{\mathbf{s}}-1}}\left(\mathfrak{P}^{(2)}_{\bb{\mathbf{s}}-1}\right)$}};
\node[] (d3b) at (6.4,-3) [] {\rotatebox{-90}{$\,=$}}; 

\node[] (e0) at (1.15,-4) [] {$\sigma \Biggl($};
\node[] (e1) at (2.75,-4) [] {$\vdots$};
\node[] (e2) at (4.5,-4) [] {$,\cdots ,$}; 
\node[] (e3) at (6.4,-4) [] {$\mathrm{sc}^{([1],2)}_{s_{\bb{\mathbf{s}}-1}}\left(\mathfrak{P}^{(2)}_{\bb{\mathbf{s}}-1}\right)$};
\node[] (e4) at (8,-4) [] {$\Biggr)$};

\node[] (f1) at (2.75,-5) [] {\rotatebox{-90}{$\,=$}};  
\node[] (f3) at (6.8,-4.95) [] {\rotatebox{-90}{
$\xymatrix@C=30pt{
\,
\ar@{=>}[r]^-{\text{
\rotatebox{90}{$\mathfrak{P}^{(2)}_{\bb{\mathbf{s}}-1}$}
}}
&
\,
}$}}; 

\node[] (g0) at (1.15,-6) [] {$\sigma \Biggl($};
\node[] (g1) at (2.75,-6) [] {$\mathrm{sc}^{([1],2)}_{s_{0}}\left(\mathfrak{Q}^{(2)}_{0}\right)$}; 
\node[] (g2) at (4.5,-6) [] {$,\cdots ,$}; 
\node[] (g3) at (6.4,-6) [] {$\mathrm{tg}^{([1],2)}_{s_{\bb{\mathbf{s}}-1}}\left(\mathfrak{P}^{(2)}_{\bb{\mathbf{s}}-1}\right)$};
\node[] (g4) at (8,-6) [] {$\Biggr)$};

\node[] (h1) at (3.05,-6.95) [] {\rotatebox{-90}{
$\xymatrix@C=30pt{
\,
\ar@{=>}[r]^-{\text{
\rotatebox{90}{$\mathfrak{Q}^{(2)}_{0}$}
}}
&
\,
}$}}; 
\node[] (h3) at (6.4,-7) [] {\rotatebox{-90}{$\,=$}};

\node[] (i0) at (1.15,-8) [] {$\sigma \Biggl($};
\node[] (i1) at (2.75,-8) [] {$\mathrm{tg}^{([1],2)}_{s_{0}}\left(\mathfrak{Q}^{(2)}_{0}\right)$}; 
\node[] (i2) at (4.5,-8) [] {$,\cdots ,$}; 
\node[] (i3) at (6.4,-8) [] {$\vdots$};
\node[] (i4) at (8,-8) [] {$\Biggr)$};

\node[] (j1) at (2.75,-9) [] {\color{white}{$\mathrm{tg}^{([1],2)}_{s_{0}}\left(\mathfrak{Q}^{(2)}_{0}\right)$}}; 
\node[] (j1b) at (2.75,-9) [] {\rotatebox{-90}{$\,=$}}; 
\node[] (j2) at (4.5,-9) [] {$\ddots$}; 
\node[] (j3) at (6.4,-9) [] {\color{white}{$\mathrm{sc}^{([1],2)}_{s_{\bb{\mathbf{s}}-1}}\left(\mathfrak{P}^{(2)}_{\bb{\mathbf{s}}-1}\right)$}};
\node[] (j3b) at (6.4,-9) [] {\rotatebox{-90}{$\,=$}}; 

\node[] (k0) at (1.15,-10) [] {$\sigma \Biggl($};
\node[] (k1) at (2.75,-10) [] {$\vdots$};
\node[] (k2) at (4.5,-10) [] {$,\cdots ,$}; 
\node[] (k3) at (6.4,-10) [] {$\mathrm{sc}^{([1],2)}_{s_{\bb{\mathbf{s}}-1}}\left(\mathfrak{P}^{(2)}_{\bb{\mathbf{s}}-1}\right)$};
\node[] (k4) at (8,-10) [] {$\Biggr)$};

\node[] (l1) at (2.75,-11) [] {\rotatebox{-90}{$\,=$}};  
\node[] (l3) at (6.8,-10.95) [] {\rotatebox{-90}{
$\xymatrix@C=30pt{
\,
\ar@{=>}[r]^-{\text{
\rotatebox{90}{$\mathfrak{Q}^{(2)}_{\bb{\mathbf{s}}-1}$}
}}
&
\,
}$}};

\node[] (m0) at (1.15,-12) [] {$\sigma \Biggl($};
\node[] (m1) at (2.75,-12) [] {$\mathrm{tg}^{([1],2)}_{s_{0}}\left(\mathfrak{Q}^{(2)}_{0}\right)$}; 
\node[] (m2) at (4.5,-12) [] {$,\cdots ,$}; 
\node[] (m3) at (6.4,-12) [] {$\mathrm{tg}^{([1],2)}_{s_{\bb{\mathbf{s}}-1}}\left(\mathfrak{P}^{(2)}_{\bb{\mathbf{s}}-1}\right)$};
\node[] (m4) at (8,-12) [] {$\Biggr)$};

\tikzset{encercla/.style={draw=black, line width=.5pt, inner sep=0pt, rectangle, rounded corners}};
\node [encercla,  fit=(a1) ] {} ;
\node [encercla,  fit=(c1)(g1) ] {} ;
\node [encercla,  fit=(i1)(m1) ] {} ;
\node [encercla,  fit=(a3)(e3) ] {} ;
\node [encercla,  fit=(g3)(k3) ] {} ;
\node [encercla,  fit=(m3) ] {} ;
\tikzset{encercla2/.style={draw=white, dashed, line width=.65pt, inner sep=0pt, rectangle, rounded corners}};
\node [encercla2,  fit=(d1) ] {} ;
\node [encercla2,  fit=(j1) ] {} ;
\node [encercla2,  fit=(d3) ] {} ;
\node [encercla2,  fit=(j3) ] {} ;
\end{tikzpicture}
\caption{
$\sigma^{\mathbf{Pth}_{\boldsymbol{\mathcal{A}}^{(2)}}}(
(
\mathfrak{Q}^{(2)}_{j}
)_{j\in \bb{\mathbf{s}}}
)
\circ^{1\mathbf{Pth}_{\boldsymbol{\mathcal{A}}^{(2)}}}
\sigma^{\mathbf{Pth}_{\boldsymbol{\mathcal{A}}^{(2)}}}(
(
\mathfrak{P}^{(2)}_{j}
)_{j\in \bb{\mathbf{s}}}
)$.
}\label{FDVVarB8Rgt}
\end{figure}

All in all, we conclude that the above second-order paths are in the kernel of the second-order Curry-Howard mapping, i.e.,
\allowdisplaybreaks
\begin{multline*}
\left[
\sigma^{\mathbf{Pth}_{\boldsymbol{\mathcal{A}}^{(2)}} }\left(
\left(
\mathfrak{Q}^{(2)}_{j}
\circ^{1\mathbf{Pth}_{\boldsymbol{\mathcal{A}}^{(2)}} }_{s_{j}}
\mathfrak{P}^{(2)}_{j}
\right)_{j\in\bb{\mathbf{s}}}
\right)
\right]_{s}
\\=
\left[
\sigma^{\mathbf{Pth}_{\boldsymbol{\mathcal{A}}^{(2)}}}\left(
\left(
\mathfrak{Q}^{(2)}_{j}
\right)_{j\in \bb{\mathbf{s}}}
\right)
\circ^{1\mathbf{Pth}_{\boldsymbol{\mathcal{A}}^{(2)}}}
\sigma^{\mathbf{Pth}_{\boldsymbol{\mathcal{A}}^{(2)}}}\left(
\left(
\mathfrak{P}^{(2)}_{j}
\right)_{j\in \bb{\mathbf{s}}}
\right)
\right]_{s}.
\end{multline*}

Therefore, their $\llbracket\cdot \rrbracket$-class must coincide, i.e.,
\allowdisplaybreaks
\begin{multline*}
\biggl\llbracket
\sigma^{\mathbf{Pth}_{\boldsymbol{\mathcal{A}}^{(2)}} }\left(
\left(
\mathfrak{Q}^{(2)}_{j}
\circ^{1\mathbf{Pth}_{\boldsymbol{\mathcal{A}}^{(2)}} }_{s_{j}}
\mathfrak{P}^{(2)}_{j}
\right)_{j\in\bb{\mathbf{s}}}
\right)
\biggr\rrbracket_{s}
\\=
\biggl\llbracket
\sigma^{\mathbf{Pth}_{\boldsymbol{\mathcal{A}}^{(2)}}}\left(
\left(
\mathfrak{Q}^{(2)}_{j}
\right)_{j\in \bb{\mathbf{s}}}
\right)
\circ^{1\mathbf{Pth}_{\boldsymbol{\mathcal{A}}^{(2)}}}
\sigma^{\mathbf{Pth}_{\boldsymbol{\mathcal{A}}^{(2)}}}\left(
\left(
\mathfrak{P}^{(2)}_{j}
\right)_{j\in \bb{\mathbf{s}}}
\right)
\biggr\rrbracket_{s}.
\end{multline*}

This completes the proof.
\end{proof}

\begin{restatable}{proposition}{PDVCtyAlg}
\label{PDVCtyAlg} $\llbracket\mathsf{Pth}_{\boldsymbol{\mathcal{A}}^{(2)}}\rrbracket$ is a $2$-categorial $\Sigma$-algebra.
\end{restatable}
\begin{proof}
That $\llbracket\mathsf{Pth}_{\boldsymbol{\mathcal{A}}^{(2)}}\rrbracket$ is an $S$-sorted $2$-category was already proven in Proposition~\ref{PDVCat}. We are only left to prove that, for every $(\mathbf{s},s)$ in $S^{\star}\times S$, $\sigma^{\llbracket\mathsf{Pth}_{\boldsymbol{\mathcal{A}}^{(2)}}\rrbracket}$ is a $2$-functor from $\llbracket\mathsf{Pth}_{\boldsymbol{\mathcal{A}}^{(2)}}\rrbracket_{\mathbf{s}}$ to $\llbracket\mathsf{Pth}_{\boldsymbol{\mathcal{A}}^{(2)}}\rrbracket_{s}$. This is a fact that follows from Propositions~\ref{PDVVarA7},~\ref{PDVVarA8},~\ref{PDVVarB7} and~\ref{PDVVarB8}.

This completes the proof.
\end{proof}

\section{
\texorpdfstring
{An Artinian preorder on $\coprod\llbracket \mathrm{Pth}_{\boldsymbol{\mathcal{A}}^{(2)}}\rrbracket$}
{An Artinian preorder on the second-order quotient}
} 

The aim of this section is to provide the coproduct of the many-sorted set of second-order path classes $\coprod \llbracket\mathrm{Pth}_{\boldsymbol{\mathcal{A}}^{(2)}}\rrbracket$ with a structure of Artinian preorder. 

\begin{restatable}{definition}{DDVOrd}
\label{DDVOrd} 
\index{partial preorder!second-order!$\leq_{\llbracket\mathbf{Pth}_{\boldsymbol{\mathcal{A}}^{(2)}}\rrbracket}$}
Let $\leq_{\llbracket\mathbf{Pth}_{\boldsymbol{\mathcal{A}}^{(2)}}\rrbracket}$ be the binary relation defined on $\coprod \llbracket\mathrm{Pth}_{\boldsymbol{\mathcal{A}}^{(2)}}\rrbracket$ containing every pair $((\llbracket\mathfrak{Q}^{(2)}\rrbracket_{t},t), (\llbracket\mathfrak{P}^{(2)}\rrbracket_{s},s))$ in $\coprod \llbracket\mathrm{Pth}_{\boldsymbol{\mathcal{A}}^{(2)}}\rrbracket$ satisfying that there exists a natural number $m\in\mathbb{N}-\{0\}$, a word $\mathbf{w}\in S^{\star}$ of length $\bb{\mathbf{w}}=m+1$ and a family of second-order paths $(\mathfrak{R}^{(2)}_{k})_{k\in\bb{\mathbf{w}}}$ in $\mathrm{Pth}_{\boldsymbol{\mathcal{A}}^{(2)},\mathbf{w}}$ such that $w_{0}=t$, $\llbracket \mathfrak{R}^{(2)}_{0}\rrbracket_{t}=\llbracket\mathfrak{Q}^{(2)}\rrbracket_{t}$, $w_{m}=s$, $\llbracket\mathfrak{R}^{(2)}_{m}\rrbracket_{s}=\llbracket\mathfrak{P}^{(2)}\rrbracket_{s}$ and, for every $k\in m$, $w_{k}=w_{k+1}$ and $\llbracket\mathfrak{R}^{(2)}_{k}\rrbracket_{w_{k}}=\llbracket\mathfrak{R}^{(2)}_{k+1}\rrbracket_{w_{m+1}}$ or $(\mathfrak{R}^{(2)}_{k}, w_{k})\leq_{\mathbf{Pth}_{\boldsymbol{\mathcal{A}}^{(2)}}} (\mathfrak{R}^{(2)}_{k+1}, w_{k+1})$.
\end{restatable}

The aim of this section is to prove that $\leq_{\llbracket\mathbf{Pth}_{\boldsymbol{\mathcal{A}}^{(2)}}\rrbracket}$ is an Artinian preorder. This is what we do in the following Proposition.

\begin{restatable}{proposition}{PDVOrdArt}
\label{PDVOrdArt} $(\coprod \llbracket\mathrm{Pth}_{\boldsymbol{\mathcal{A}}^{(2)}}\rrbracket, \leq_{\llbracket\mathbf{Pth}_{\boldsymbol{\mathcal{A}}^{(2)}}\rrbracket})$ is a preordered set. Moreover, in this preordered set there is not any strictly decreasing $\omega_{0}$-chain, i.e., $(\coprod \llbracket \mathrm{Pth}_{\boldsymbol{\mathcal{A}}^{(2)}}\rrbracket, \leq_{\llbracket \mathbf{Pth}_{\boldsymbol{\mathcal{A}}^{(2)}}\rrbracket})$ is an Artinian partially preordered set.
\end{restatable}
\begin{proof}

That $\leq_{\llbracket\mathbf{Pth}_{\boldsymbol{\mathcal{A}}^{(2)}}\rrbracket}$ is reflexive follows from the fact that either $\llbracket \cdot \rrbracket$ or $\leq_{\mathbf{Pth}_{\boldsymbol{\mathcal{A}}^{(2)}}}$ are reflexive. 

Indeed, let $s$ be a sort in $S$ and let $\llbracket\mathfrak{P}^{(2)}\rrbracket_{s}$ be a second-order path class in $\llbracket\mathrm{Pth}_{\boldsymbol{\mathcal{A}}^{(2)}}\rrbracket_{s}$.  Then, for the sequence $(\mathfrak{P}^{(2)},\mathfrak{P}^{(2)})$ in $\mathrm{Pth}_{\boldsymbol{\mathcal{A}}^{(2)},ss}$, we have that

\[
\left(
\left\llbracket \mathfrak{P}^{(2)}\right\rrbracket_{s}, s\right)
\leq_{\llbracket \mathbf{Pth}_{\boldsymbol{\mathcal{A}}^{(2)}}\rrbracket} 
\left(
\left\llbracket \mathfrak{P}^{(2)}\right\rrbracket_{s}, s\right).
\]

That $\leq_{\llbracket \mathbf{Pth}_{\boldsymbol{\mathcal{A}}^{(2)}}\rrbracket}$ is transitive follows as well. Indeed, let $s,t,u$ be sorts in $S$ and let $\llbracket\mathfrak{P}^{(2)}\rrbracket_{s}$, $\llbracket \mathfrak{Q}^{(2)}\rrbracket_{t}$ and $\llbracket\mathfrak{S}^{(2)}\rrbracket_{u}$ be second-order path classes satisfying that 
\[
\left(
\left\llbracket 
\mathfrak{Q}^{(2)}
\right\rrbracket_{t},t\right)
\leq_{\llbracket \mathbf{Pth}_{\boldsymbol{\mathcal{A}}^{(2)}}\rrbracket} 
\left(
\left\llbracket 
\mathfrak{P}^{(2)}
\right\rrbracket_{s},s\right)
\]
\[
\left(
\left\llbracket 
\mathfrak{S}^{(2)}
\right\rrbracket_{u},u\right)
\leq_{\llbracket \mathbf{Pth}_{\boldsymbol{\mathcal{A}}^{(2)}}\rrbracket} 
\left(
\left\llbracket 
\mathfrak{Q}^{(2)}
\right\rrbracket_{t},t\right)
\]
Thus there exists a natural number $m\in\mathbb{N}-\{0\}$, a word $\mathbf{w}\in S^{\star}$ of length $\bb{\mathbf{w}}=m+1$ and a family of second-order paths $(\mathfrak{R}^{(2)}_{k})_{k\in\bb{\mathbf{w}}}$ in $\mathrm{Pth}_{\boldsymbol{\mathcal{A}}^{(2)},\mathbf{w}}$ such that $w_{0}=t$, $\llbracket\mathfrak{R}^{(2)}_{0}\rrbracket_{t}=\llbracket \mathfrak{Q}^{(2)}\rrbracket_{t}$, $w_{m}=s$, $\llbracket \mathfrak{R}^{(2)}_{m}\rrbracket_{s}=\llbracket \mathfrak{P}^{(2)}\rrbracket_{s}$ and, for every $k\in m$, $w_{k}=w_{k+1}$ and $\llbracket \mathfrak{R}^{(2)}_{k}\rrbracket_{w_{k}}= \llbracket \mathfrak{R}^{(2)}_{k+1}\rrbracket_{w_{k+1}}$ or $(\mathfrak{R}^{(2)}_{k}, w_{k})\leq_{\mathbf{Pth}_{\boldsymbol{\mathcal{A}}^{(2)}}} (\mathfrak{R}^{(2)}_{k+1}, w_{k+1})$. Regarding the other inequality, there exists a natural number $m'\in\mathbb{N}-\{0\}$, a word $\mathbf{w}'\in S^{\star}$ of length $\bb{\mathbf{w}'}=m'+1$ and a family of second-order paths $(\mathfrak{R}'^{(2)}_{k})_{k\in\bb{\mathbf{w}}}$ in $\mathrm{Pth}_{\boldsymbol{\mathcal{A}}^{(2)},\mathbf{w}'}$ such that $w'_{0}=u$, $\llbracket\mathfrak{R}'^{(2)}_{0}\rrbracket_{u}=\llbracket \mathfrak{S}^{(2)}\rrbracket_{u}$, $w_{m'}=t$, $\llbracket \mathfrak{R}'^{(2)}_{m'}\rrbracket_{t}=\llbracket \mathfrak{Q}^{(2)}\rrbracket_{t}$ and, for every $k\in m'$, $w'_{k}=w'_{k+1}$ and $\llbracket \mathfrak{R}'^{(2)}_{k}\rrbracket_{w'_{k}}= \llbracket \mathfrak{R}'^{(2)}_{k+1}\rrbracket_{w'_{k+1}}$ or $(\mathfrak{R}'^{(2)}_{k}, w'_{k})\leq_{\mathbf{Pth}_{\boldsymbol{\mathcal{A}}^{(2)}}} (\mathfrak{R}'^{(2)}_{k+1}, w'_{k+1})$

Thus, the concatenation of the sequences $(\mathfrak{R}^{(2)}_{k})_{k\in\bb{\mathbf{w}}}$ and $(\mathfrak{R}'^{(2)}_{k})_{k\in\bb{\mathbf{w}'}}$ instantiates that 
\[
\left(
\left\llbracket 
\mathfrak{S}^{(2)}
\right\rrbracket_{u},u\right)
\leq_{\llbracket \mathbf{Pth}_{\boldsymbol{\mathcal{A}}^{(2)}}\rrbracket} 
\left(
\left\llbracket 
\mathfrak{P}^{(2)}
\right\rrbracket_{s},s\right).
\]

This proves that 
$\leq_{\llbracket \mathbf{Pth}_{\boldsymbol{\mathcal{A}}^{(2)}}\rrbracket}$ is a preorder on  $\coprod \llbracket\mathrm{Pth}_{\boldsymbol{\mathcal{A}}^{(2)}}\rrbracket$.

That it is Artinian follows from the fact that any strictly decreasing chain in $(\coprod \llbracket\mathrm{Pth}_{\boldsymbol{\mathcal{A}}^{(2)}}\rrbracket, \leq_{\llbracket\mathbf{Pth}_{\boldsymbol{\mathcal{A}}^{(2)}}\rrbracket})$ is given by a strictly decreasing chain that only uses the strict order $<_{\mathbf{Pth}_{\boldsymbol{\mathcal{A}}^{(2)}}}$ at each step. The statement follows from the fact that, by Proposition~\ref{PDOrdArt}, $(\coprod\mathrm{Pth}_{\boldsymbol{\mathcal{A}}^{(2)}}, \leq_{\mathbf{Pth}_{\boldsymbol{\mathcal{A}}^{(2)}}})$ is an Artinian preordered set.

This completes the proof.
\end{proof}

Taking into account the definitions above, we now prove that the coproduct of the projection mapping becomes order-preserving.

\begin{restatable}{proposition}{PDVProjMono}
\label{PDVProjMono} The mapping $\coprod\mathrm{pr}^{\llbracket\cdot\rrbracket }$ is monotone
\[
\textstyle
\coprod\mathrm{pr}^{\llbracket\cdot\rrbracket}\colon 
\left(
\coprod \mathrm{Pth}_{\boldsymbol{\mathcal{A}}^{(2)}},
\leq_{\mathbf{Pth}_{\boldsymbol{\mathcal{A}}^{(2)}}}
\right)
\mor
\left(
\coprod \llbracket \mathrm{Pth}_{\boldsymbol{\mathcal{A}}^{(2)}}\rrbracket,
\leq_{\llbracket\mathbf{Pth}_{\boldsymbol{\mathcal{A}}^{(2)}}\rrbracket}
\right)
\]
\end{restatable}
\begin{proof}
Let $s$ and $t$ be sorts in $S$ and let $(\mathfrak{Q}^{(2)},t)$ and $(\mathfrak{P}^{(2)},s)$ be pairs in $\coprod \mathrm{Pth}_{\boldsymbol{\mathcal{A}}^{(2)}}$ satisfying that 
\[
\left(
\mathfrak{Q}^{(2)},t
\right)
\leq_{\mathbf{Pth}_{\boldsymbol{\mathcal{A}}^{(2)}}}
\left(
\mathfrak{P}^{(2)},s
\right),
\]
then the sequence $(\mathfrak{Q}^{(2)},\mathfrak{P}^{(2)})$ can be used to instantiate that 
\[
\left(
\llbracket\mathfrak{Q}^{(2)}\rrbracket_{t},t
\right)
\leq_{\llbracket \mathbf{Pth}_{\boldsymbol{\mathcal{A}}^{(2)}}\rrbracket}
\left(
\llbracket\mathfrak{P}^{(2)}\rrbracket_{s},s
\right).
\]

It follows that the mapping $\coprod\mathrm{pr}^{\llbracket\cdot \rrbracket}$ is monotone.
\end{proof}

It also happens that the coproduct of the intermediate projections are also order-preserving.

\begin{proposition}\label{PDCHVProjMono} The mapping $\coprod\mathrm{pr}^{\llbracket\cdot\rrbracket,[\cdot]}$ is monotone
\[
\textstyle
\coprod\mathrm{pr}^{\llbracket\cdot\rrbracket, [\cdot]}\colon 
\left(
\coprod 
\left[
\mathrm{Pth}_{\boldsymbol{\mathcal{A}}^{(2)}}
\right],
\leq_{\left[\mathbf{Pth}_{\boldsymbol{\mathcal{A}}^{(2)}}\right]}
\right)
\mor
\left(
\coprod \llbracket \mathrm{Pth}_{\boldsymbol{\mathcal{A}}^{(2)}}\rrbracket,
\leq_{\llbracket\mathbf{Pth}_{\boldsymbol{\mathcal{A}}^{(2)}}\rrbracket}
\right)
\]
\end{proposition}
\begin{proof}
Let $s$ and $t$ be sorts in $S$ and let $([\mathfrak{Q}^{(2)}]_{t},t)$ and $([\mathfrak{P}^{(2)}]_{s},s)$ be pairs in $\coprod [\mathrm{Pth}_{\boldsymbol{\mathcal{A}}^{(2)}}]$ satisfying that 
\[
\left(
\left[\mathfrak{Q}^{(2)}\right]_{t},t
\right)
\leq_{[\mathbf{Pth}_{\boldsymbol{\mathcal{A}}^{(2)}}]}
\left(
\left[\mathfrak{P}^{(2)}\right]_{s},s
\right),
\]
then, according to Definition~\ref{DDCHOrd},
there exists second-order path representatives $\mathfrak{Q}'^{(2)}\in [\mathfrak{Q}^{(2)}]_{t}$ and $\mathfrak{P}'^{(2)}\in [\mathfrak{P}^{(2)}]_{s}$ satisfying that 
\[
\left(
\mathfrak{Q}'^{(2)},t
\right)
\leq_{\mathbf{Pth}_{\boldsymbol{\mathcal{A}}^{(2)}}}
\left(
\mathfrak{P}^{(2)},s
\right).
\]

Thus, the sequence $(\mathfrak{Q}'^{(2)},\mathfrak{P}'^{(2)})$ can be used to instantiate that 
\[
\left(
\llbracket\mathfrak{Q}^{(2)}\rrbracket_{t},t
\right)
\leq_{\llbracket \mathbf{Pth}_{\boldsymbol{\mathcal{A}}^{(2)}}\rrbracket}
\left(
\llbracket\mathfrak{P}^{(2)}\rrbracket_{s},s
\right).
\]

It follows that the mapping $\coprod\mathrm{pr}^{\llbracket\cdot \rrbracket, [\cdot]}$ is monotone.
\end{proof}

\begin{proposition}\label{PDUpsVProjMono} The mapping $\coprod\mathrm{pr}^{\llbracket\cdot\rrbracket,\Upsilon^{[1]}}$ is monotone
\[
\textstyle
\coprod\mathrm{pr}^{\llbracket\cdot\rrbracket, \Upsilon^{[1]}}\colon 
\left(
\coprod 
\left[
\mathrm{Pth}_{\boldsymbol{\mathcal{A}}^{(2)}}
\right]_{\Upsilon^{[1]}},
\leq_{\left[\mathbf{Pth}_{\boldsymbol{\mathcal{A}}^{(2)}}\right]_{\Upsilon^{[1]}}}
\right)
\mor
\left(
\coprod \llbracket \mathrm{Pth}_{\boldsymbol{\mathcal{A}}^{(2)}}\rrbracket,
\leq_{\llbracket\mathbf{Pth}_{\boldsymbol{\mathcal{A}}^{(2)}}\rrbracket}
\right)
\]
\end{proposition}
\begin{proof}
Let $s$ and $t$ be sorts in $S$ and let $([\mathfrak{Q}^{(2)}]_{\Upsilon^{[1]}_{t}},t)$ and $([\mathfrak{P}^{(2)}]_{\Upsilon^{[1]}_{s}},s)$ be pairs in $\coprod [\mathrm{Pth}_{\boldsymbol{\mathcal{A}}^{(2)}}]_{\Upsilon^{[1]}}$ satisfying that 
\[
\left(
\left[\mathfrak{Q}^{(2)}\right]_{\Upsilon^{[1]}_{t}},t
\right)
\leq_{[\mathbf{Pth}_{\boldsymbol{\mathcal{A}}^{(2)}}]_{\Upsilon^{[1]}}}
\left(
\left[\mathfrak{P}^{(2)}\right]_{\Upsilon^{[1]}_{t}},s
\right),
\]
then, according to Definition~\ref{DDUpsOrd},
there exists a natural number $m\in\mathbb{N}-\{0\}$, a word $\mathbf{w}\in S^{\star}$ of length $\bb{\mathbf{w}}=m+1$ and a family of second-order paths $(\mathfrak{R}^{(2)}_{k})_{k\in\bb{\mathbf{w}}}$ in $\mathrm{Pth}_{\boldsymbol{\mathcal{A}}^{(2)},\mathbf{w}}$ such that $w_{0}=t$, $[\mathfrak{R}^{(2)}_{0}]_{\Upsilon^{[1]}_{t}}=[\mathfrak{Q}^{(2)}]_{\Upsilon^{[1]}_{t}}$, $w_{m}=s$, $[\mathfrak{R}^{(2)}_{m}]_{\Upsilon^{[1]}_{s}}=[\mathfrak{P}^{(2)}]_{\Upsilon^{[1]}_{s}}$ and, for every $k\in m$, $w_{k}=w_{k+1}$ and $(\mathfrak{R}^{(2)}_{k}, \mathfrak{R}^{(2)}_{k+1})\in \Upsilon^{[1]}_{w_{k}}$ or $(\mathfrak{R}^{(2)}_{k}, w_{k})\leq_{\mathbf{Pth}_{\boldsymbol{\mathcal{A}}^{(2)}}} (\mathfrak{R}^{(2)}_{k+1}, w_{k+1})$. Thus, the same sequence $(\mathfrak{R}^{(2)}_{k})_{k\in\bb{\mathbf{w}}}$ can be used to instantiate that 
\[
\left(
\llbracket\mathfrak{Q}^{(2)}\rrbracket_{t},t
\right)
\leq_{\llbracket \mathbf{Pth}_{\boldsymbol{\mathcal{A}}^{(2)}}\rrbracket}
\left(
\llbracket\mathfrak{P}^{(2)}\rrbracket_{s},s
\right).
\]

It follows that the mapping $\coprod\mathrm{pr}^{\llbracket\cdot \rrbracket, \Upsilon^{[1]}}$ is monotone.
\end{proof}

\chapter{
\texorpdfstring
{The free completion of the second-order identity path mapping}
{The second-order free completion}
}\label{S2I}

In this chapter we introduce the free completion of the many-sorted partial $\Sigma^{\boldsymbol{\mathcal{A}}^{(2)}}$-algebra $\mathbf{Pth}_{\boldsymbol{\mathcal{A}}^{(2)}}$, denoted by $\mathbf{F}_{\Sigma^{\boldsymbol{\mathcal{A}}^{(2)}}}(\mathbf{Pth}_{\boldsymbol{\mathcal{A}}^{(2)}})$. This allows us to consider the free completion of the $\mathrm{ip}^{(2,X)}$ mapping, denoted by $\mathrm{ip}^{(2,X)@}$, which is, by construction, a $\Sigma^{\boldsymbol{\mathcal{A}}^{(2)}}$-homomorphism. By using this mapping we show that, for each sort $s\in S$ and for each second-order path $\mathfrak{P}^{(2)}\in \mathrm{Pth}_{\boldsymbol{\mathcal{A}}^{(2)},s}$, the object $\mathrm{ip}^{(2,X)@}_{s}(\mathrm{CH}^{(2)}_{s}(\mathfrak{P}^{(2)}))$ which, a priori, is an object in  $\mathbf{F}_{\Sigma^{\boldsymbol{\mathcal{A}}^{(2)}}}(\mathbf{Pth}_{\boldsymbol{\mathcal{A}}^{(2)}})_{s}$ is indeed a proper second-order path in $\mathrm{Pth}_{\boldsymbol{\mathcal{A}}^{(2)},s}$. Moreover, we show that the second-order path $\mathrm{ip}^{(2,X)@}_{s}(\mathrm{CH}^{(2)}_{s}(\mathfrak{P}^{(2)}))$ is, in fact, a second-order path in $[\mathfrak{P}^{(2)}]_{s}$. This has a number of consequences. We also show that $\mathrm{ip}^{(2,X)@}_{s}(\mathrm{CH}^{(2)}_{s}(\mathfrak{P}^{(2)}))$ is equal to $\mathfrak{P}^{(2)}$, for a $(2,[1])$-identity second-order path $\mathfrak{P}^{(2)}$.  We prove that the composition $\mathrm{ip}^{(2,X)@}\circ \mathrm{CH}^{(2)}$ is a $\Sigma$-homomorphism that also preserves rewrite rules, the $0$-composition, the $0$-source and the $0$-target operation symbols. It also preserves second-order rewrite rules, the $1$-source and the $1$-target operation symbols.  Furthermore, for a one-step second-order path $\mathfrak{P}^{(2)}$ in $\mathrm{Pth}_{\boldsymbol{\mathcal{A}}^{(2)},s}$, we also obtain that $\mathrm{ip}^{(2,X)@}_{s}(\mathrm{CH}^{(2)}_{s}(\mathfrak{P}^{(2)}))$ is equal to $\mathfrak{P}^{(2)}$. We also prove that if $\mathfrak{P}^{(2)}$ is a $(2,[1])$-identity second-order path in $\mathrm{Pth}_{\boldsymbol{\mathcal{A}}^{(2)},s}$, then the second-order paths $\mathfrak{P}^{(2)}$ and $\mathrm{ip}^{(2,X)@}_{s}(\mathrm{CH}^{(2)}_{s}(\mathfrak{P}^{(2)}))$ are $\Upsilon^{[1]}_{s}$-related. It is also shown that, for a pair of second-order paths $\mathfrak{P}^{(2)}$ and $\mathfrak{Q}^{(2)}$ in $\mathrm{Pth}_{\boldsymbol{\mathcal{A}}^{(2)},s}$, if they are $\Upsilon^{[1]}_{s}$-related, then the second-order paths $\mathrm{ip}^{(2,X)@}_{s}(\mathrm{CH}^{(2)}_{s}(\mathfrak{P}^{(2)}))$ and $\mathrm{ip}^{(2,X)@}_{s}(\mathrm{CH}^{(2)}_{s}(\mathfrak{Q}^{(2)}))$ are also $\Upsilon^{[1]}_{s}$-related. We conclude this chapter by showing that the mappings $\mathrm{ip}^{(2,X)@}$  and  $\mathrm{ip}^{(2,X)@}\circ \mathrm{CH}^{(2)}$ are order-preserving.


We first apply this construction to the many-sorted partial $\Sigma^{\boldsymbol{\mathcal{A}}^{(2)}}$-algebra of second-order paths and we will stablish the basic results for understanding its behaviour. The following definition and the statements in it follow directly from Proposition~\ref{PFreeComp}.

\begin{restatable}{definition}{DDFP}
\label{DDFP} 
\index{free completion!second-order!$\mathbf{F}_{\Sigma^{\boldsymbol{\mathcal{A}}^{(2)}}}(\mathbf{Pth}_{\boldsymbol{\mathcal{A}}^{(2)}})$}
We let $\mathbf{F}_{\Sigma^{\boldsymbol{\mathcal{A}}^{(2)}}}(\mathbf{Pth}_{\boldsymbol{\mathcal{A}}^{(2)}})$ stand for the free $\Sigma^{\boldsymbol{\mathcal{A}}^{(2)}}$-completion of the partial $\Sigma^{\boldsymbol{\mathcal{A}}^{(2)}}$-algebra $\mathbf{Pth}_{\boldsymbol{\mathcal{A}}^{(2)}}$, and 
\index{inclusion!second-order!$\eta^{(2,\mathbf{Pth}_{\boldsymbol{\mathcal{A}}^{(2)}})}$}
we denote by $\eta^{(2,\mathbf{Pth}_{\boldsymbol{\mathcal{A}}^{(2)}})}$ the $\Sigma^{\boldsymbol{\mathcal{A}}^{(2)}}$-homomorphism
$$
\eta^{(2,\mathbf{Pth}_{\boldsymbol{\mathcal{A}}^{(2)}})}\colon
\mathbf{Pth}_{\boldsymbol{\mathcal{A}}^{(2)}}
\mor
\mathbf{F}_{\Sigma^{\boldsymbol{\mathcal{A}}^{(2)}}}\left(
\mathbf{Pth}_{\boldsymbol{\mathcal{A}}^{(2)}}
\right)
$$
given by the insertion of generators. Note that $\eta^{(2,\mathbf{Pth}_{\boldsymbol{\mathcal{A}}^{(2)}})}$ is an embedding. As usual, we identify $\mathrm{Pth}_{\boldsymbol{\mathcal{A}}^{(2)}}$ with its image under $\eta^{(2,\mathbf{Pth}_{\boldsymbol{\mathcal{A}}^{(2)}})}$. In this way, $\mathbf{Pth}_{\boldsymbol{\mathcal{A}}^{(2)}}$ becomes a weak subalgebra of $\mathbf{F}_{\Sigma^{\boldsymbol{\mathcal{A}}^{(2)}}}(\mathbf{Pth}_{\boldsymbol{\mathcal{A}}^{(2)}})$.
\end{restatable}

\begin{remark}\label{RDFP} 
Let us recall from Proposition~\ref{PFreeComp} that, for every $\Sigma^{\boldsymbol{\mathcal{A}}^{(2)}}$-algebra $\mathbf{B}$ and every homomorphism $f$ from $\mathbf{Pth}_{\boldsymbol{\mathcal{A}}^{(2)}}$ to $\mathbf{B}$, there exists a unique $\Sigma^{\boldsymbol{\mathcal{A}}^{(2)}}$-homomorphism $f^{\mathrm{fc}}$ from $\mathbf{F}_{\Sigma^{\boldsymbol{\mathcal{A}}^{(2)}}}(\mathbf{Pth}_{\boldsymbol{\mathcal{A}}^{(2)}})$ to $\mathbf{B}$ such that $f=f^{\mathrm{fc}}\circ\eta^{(2,\mathbf{Pth}_{\boldsymbol{\mathcal{A}}^{(2)}})}$.
\end{remark}

We apply the same construction to the discrete $\Sigma^{\boldsymbol{\mathcal{A}}^{(2)}}$-algebra on $X$, i.e., the many-sorted partial $\Sigma^{\boldsymbol{\mathcal{A}}^{(2)}}$-algebra whose underlying set is $X$ and whose operations are defined nowhere.

\begin{definition}\label{DDFX} We let $\mathbf{D}_{\Sigma^{\boldsymbol{\mathcal{A}}^{(2)}}}(X)$ stand for the partial $\Sigma^{\boldsymbol{\mathcal{A}}^{(2)}}$-algebra whose underlying $S$-sorted set is $X$ and which is such that, for each $(\mathbf{s},s)\in S^{\star}\times S$ and for each operation symbol in $\Sigma^{\boldsymbol{\mathcal{A}}^{(2)}}_{\mathbf{s},s}$, the domain of definition of the corresponding operation is empty. We call it the \emph{discrete} $\Sigma^{\boldsymbol{\mathcal{A}}^{(2)}}$-algebra \emph{on} $X$. For $\mathbf{D}_{\Sigma^{\boldsymbol{\mathcal{A}}^{(2)}}}(X)$ it happens that every $S$-sorted mapping $f$ from $X$ to a partial $\Sigma^{\boldsymbol{\mathcal{A}}^{(2)}}$-algebra $\mathbf{B}$ is a $\Sigma^{\boldsymbol{\mathcal{A}}^{(2)}}$-homomorphism, 
$$
f\colon\mathbf{D}_{\Sigma^{\boldsymbol{\mathcal{A}}^{(2)}}}(X)
\mor
\mathbf{B}.
$$
Moreover, by Remark~\ref{RFreeDisc}, $\mathbf{F}_{\Sigma^{\boldsymbol{\mathcal{A}}^{(2)}}}(
\mathbf{D}_{\Sigma^{\boldsymbol{\mathcal{A}}^{(2)}}}(X))$, the free completion of $\mathbf{D}_{\Sigma^{\boldsymbol{\mathcal{A}}^{(2)}}}(X)$, and $\mathbf{T}_{\Sigma^{\boldsymbol{\mathcal{A}}^{(2)}}}(X)$, the free $\Sigma^{\boldsymbol{\mathcal{A}}^{(2)}}$-algebra on $X$, are such that $\mathbf{F}_{\Sigma^{\boldsymbol{\mathcal{A}}^{(2)}}}(
\mathbf{D}_{\Sigma^{\boldsymbol{\mathcal{A}}^{(2)}}}(X))= \mathbf{T}_{\Sigma^{\boldsymbol{\mathcal{A}}^{(2)}}}(X)$. For this reason, the standard insertion of generators $\eta^{(2,X)}$ from $X$ to $\mathrm{T}_{\Sigma^{\boldsymbol{\mathcal{A}}^{(2)}}}(X)$ becomes a $\Sigma^{\boldsymbol{\mathcal{A}}^{(2)}}$-homomorphism
$$
\eta^{(2,X)}\colon
\mathbf{D}_{\Sigma^{\boldsymbol{\mathcal{A}}^{(2)}}}(X)
\mor
\mathbf{T}_{\Sigma^{\boldsymbol{\mathcal{A}}^{(2)}}}(X).
$$
Let us note that $\eta^{(2,X)}$ is an embedding. As usual, we identify $\mathbf{D}_{\Sigma^{\boldsymbol{\mathcal{A}}^{(2)}}}(X)$ with its image under $\eta^{(2,X)}$. In this way, $\mathbf{D}_{\Sigma^{\boldsymbol{\mathcal{A}}^{(2)}}}(X)$ becomes a weak subalgebra of $\mathbf{T}_{\Sigma^{\boldsymbol{\mathcal{A}}^{(2)}}}(X)$.
\end{definition}

\begin{remark}\label{RDFX} 
Let us recall from Proposition~\ref{PFreeComp} that, for every $\Sigma^{\boldsymbol{\mathcal{A}}^{(2)}}$-algebra $\mathbf{B}$ and every homomorphism $f$ from $\mathbf{D}_{\Sigma^{\boldsymbol{\mathcal{A}}^{(2)}}}(X)$ to $\mathbf{B}$, there exists a unique $\Sigma^{\boldsymbol{\mathcal{A}}^{(2)}}$-homomorphism $f^{\mathrm{fc}}$ from $\mathbf{F}_{\Sigma^{\boldsymbol{\mathcal{A}}^{(2)}}}(\mathbf{D}_{\Sigma^{\boldsymbol{\mathcal{A}}^{(2)}}}(X))$ to $\mathbf{B}$ such that $f=f^{\mathrm{fc}}\circ\eta^{(2,X)}$. The reader must acknowledge that this is the universal property for the free $\Sigma^{\boldsymbol{\mathcal{A}}^{(2)}}$-algebra on $X$. 
\end{remark}

We next consider the mapping $\mathrm{ip}^{(2,X)}$ from $X$ to $\mathbf{Pth}_{\boldsymbol{\mathcal{A}}^{(2)}}$, that, for a sort $s\in S$ and a variable $x\in X_{s}$ associates the $(2,[1])$-identity second-order path on $[x]_{s}$. For the discrete many-sorted partial $\Sigma^{\boldsymbol{\mathcal{A}}^{(2)}}$-algebra structure on $X$ this simple many-sorted mapping becomes a $\Sigma^{\boldsymbol{\mathcal{A}}^{(2)}}$-homomorphism. The next statements follow directly from Propositions~\ref{PFreeComp} and~\ref{PSch}.

\begin{restatable}{definition}{DDIp}
\label{DDIp} 
The $S$-sorted mapping $\mathrm{ip}^{(2,X)}$ from $X$ to $\mathbf{Pth}_{\boldsymbol{\mathcal{A}}}$ that, for every  sort $s\in S$, sends a variable $x\in X_{s}$ to $\mathrm{ip}^{(2,[1])\sharp}_{s}([x]_{s})$, the $(2,[1])$-identity second-order path on $[x]_{s}$ is a $\Sigma^{\boldsymbol{\mathcal{A}}^{(2)}}$-homomorphism from  $\mathbf{D}_{\Sigma^{\boldsymbol{\mathcal{A}}^{(2)}}}(X)$
$$
\mathrm{ip}^{(2,X)}
\colon
\mathbf{D}_{\Sigma^{\boldsymbol{\mathcal{A}}^{(2)}}}(X)
\mor
\mathbf{Pth}_{\boldsymbol{\mathcal{A}}^{(2)}}.
$$
\index{identity!second-order!$\mathrm{ip}^{(2,X)\mathrm{fc}}$}
On the other hand, by Corollary~\ref{CFreeAdj}, there exists a unique $\Sigma^{\boldsymbol{\mathcal{A}}^{(2)}}$-homomorphism, denoted $\mathrm{ip}^{(2,X)@}$, from $\mathbf{T}_{\Sigma^{\boldsymbol{\mathcal{A}}^{(2)}}}(X)$ to $\mathbf{F}_{\Sigma^{\boldsymbol{\mathcal{A}}^{(2)}}}(\mathbf{Pth}_{\boldsymbol{\mathcal{A}}^{(2)}})$,
$$
\mathrm{ip}^{(2,X)@}
\colon
\mathbf{T}_{\Sigma^{\boldsymbol{\mathcal{A}}^{(2)}}}(X)
\mor
\mathbf{F}_{\Sigma^{\boldsymbol{\mathcal{A}}^{(2)}}}(\mathbf{Pth}_{\boldsymbol{\mathcal{A}}^{(2)}})
$$ 
such that 
$
\mathrm{ip}^{(2,X)@}\circ
\eta^{(2,X)}
=
\eta^{(2,\mathbf{Pth}_{\boldsymbol{\mathcal{A}}^{(2)}})}
\circ
\mathrm{ip}^{(2,X)}.
$
We will refer to $\mathrm{ip}^{(2,X)@}$ as the \emph{free completion of the $S$-sorted mapping $\mathrm{ip}^{(2,X)}$}. 
In this particular case we have that 
$$
\mathrm{ip}^{(2,X)@}=\left(
\eta^{(2,\mathbf{Pth}_{\boldsymbol{\mathcal{A}}^{(2)}})}\circ \mathrm{ip}^{(2,X)}
\right)^{\mathrm{fc}}
$$
and that
$$
\mathbf{Sch}\left(\mathrm{ip}^{(2,X)}
\right)=\left(\left(\eta^{(2,\mathbf{Pth}_{\boldsymbol{\mathcal{A}}^{(2)}})}\circ \mathrm{ip}^{(2,X)}\right)^{\mathrm{fc}}
\right)^{-1}
\left[\eta^{(2,\mathbf{Pth}_{\boldsymbol{\mathcal{A}}^{(2)}})}
\left[\mathbf{Pth}_{\boldsymbol{\mathcal{A}}^{(2)}}
\right]\right].
$$
where, by Proposition~\ref{PSch}, $\mathbf{Sch}(\mathrm{ip}^{(2,X)})$ is the $X$-generated relative subalgebra of the free $\Sigma^{\boldsymbol{\mathcal{A}}^{(2)}}$-completion of $\mathbf{D}_{\Sigma^{\boldsymbol{\mathcal{A}}^{(2)}}}(X)$, which, as we know, is $\mathbf{T}_{\Sigma^{\boldsymbol{\mathcal{A}}^{(2)}}}(X)$. Finally, we let $\mathrm{ip}^{(2,X)\mathrm{Sch}}$ stand for the unique closed $\Sigma^{\boldsymbol{\mathcal{A}}^{(2)}}$-homomorphism from $\mathbf{Sch}(\mathrm{ip}^{(2,X)})$ to $\mathbf{Pth}_{\boldsymbol{\mathcal{A}}^{(2)}}$ making the square of the diagram
in Figure~\ref{FDIp} commute.
\end{restatable}

\begin{figure}
\begin{tikzpicture}
[ACliment/.style={-{To [angle'=45, length=5.75pt, width=4pt, round]}
}]
\node[] (A) 	at 	(-3,1.5) 	[] 	{$\mathbf{D}_{\Sigma^{\boldsymbol{\mathcal{A}}^{(2)}}}(X)$};
\node[] (B) 	at 	(0,-2) 	[] 	{$\mathbf{Pth}_{\boldsymbol{\mathcal{A}}^{(2)}}$};
\node[]	(S)		at 	(0,0)	[]	{$\mathbf{Sch}(\mathrm{ip}^{(2,X)})$};
\node[] (FA)	at	(4,0)	[]	{$\mathbf{T}_{\Sigma^{\boldsymbol{\mathcal{A}}^{(2)}}}(X)$};
\node[] (FB)	at	(4,-2)	[]	{$\mathbf{F}_{\Sigma^{\boldsymbol{\mathcal{A}}^{(2)}}}(\mathbf{Pth}_{\boldsymbol{\mathcal{A}}^{(2)}})$};
\draw[ACliment, bend right]  (A) 	to node [below left]	{$\mathrm{ip}^{(2,X)}$} (B);
\draw[ACliment, bend left]  (A) 	to node [above right]	{$\eta^{(2,X)}$} (FA);
\draw[ACliment]  (A) 	to node [above right]
{$\mathrm{in}^{\mathbf{D}_{\Sigma^{\boldsymbol{\mathcal{A}}^{(2)}}}(X),\mathbf{Sch}(\mathrm{ip}^{(2,X)})}$} (S);
\draw[ACliment]  (B) 	to node [below]	{$\eta^{(2,\mathbf{Pth}_{\boldsymbol{\mathcal{A}}^{(2)}})}$} (FB);
\draw[ACliment]  (S) 	to node [above]	{$\mathrm{in}^{\mathbf{Sch}(\mathrm{ip}^{(2,X)})}$} (FA);
\draw[ACliment]  (S) 	to node [right]	{$\mathrm{ip}^{(2,X)\mathrm{Sch}}$} (B);
\draw[ACliment]  (FA) 	to node [right]	{$\mathrm{ip}^{(2,X)@}$} (FB);
\end{tikzpicture}
\caption{The mapping $\mathrm{ip}^{(2,X)@}$.}
\label{FDIp}
\end{figure}


\section{
\texorpdfstring
{Interaction with the second-order Curry-Howard mapping}
{Interaction}
}
Once we have introduced the mapping $\mathrm{ip}^{(2,X)@}$ mapping, we will present a series of results aiming at having a better insight on its behaviour with the terms in $\mathrm{T}_{\Sigma^{\boldsymbol{\mathcal{A}}^{(2)}}}(X)$ that are in the image of $\mathrm{Pth}_{\boldsymbol{\mathcal{A}}^{(2)}}$ by means of the second-order Curry-Howard mapping.

We first explore the relation between the mapping $\mathrm{ip}^{(2,X)@}$ and the value of the second-order Curry-Howard mapping at  $(2,[1])$-identity second-order paths. Before doing so, we introduce the $\Sigma^{\boldsymbol{\mathcal{A}}}$-reduct of the free $\Sigma^{\boldsymbol{\mathcal{A}}^{(2)}}$-completion of the many-sorted partial $\Sigma^{\boldsymbol{\mathcal{A}}^{(2)}}$-algebra of second-order paths

\begin{definition}\label{DDPFRedDU} Let $\mathrm{in}^{\Sigma,(2,1)}$ be the canonical embedding of $\Sigma^{\boldsymbol{\mathcal{A}}}$ into $\Sigma^{\boldsymbol{\mathcal{A}}^{(2)}}$. Then by Proposition~\ref{PFunSig}, for the morphism $\mathbf{in}^{\Sigma,(2,1)}=(\mathrm{id}_{S},\mathrm{in}^{\Sigma,(2,1)})$ from $(S,\Sigma^{\boldsymbol{\mathcal{A}}})$ to $(S,\Sigma^{\boldsymbol{\mathcal{A}}^{(2)}})$ and the free $\Sigma^{\boldsymbol{\mathcal{A}}^{(2)}}$-completion $\mathbf{F}_{\Sigma^{\boldsymbol{\mathcal{A}}^{(2)}}}(\mathbf{Pth}_{\boldsymbol{\mathcal{A}}^{(2)}})$, we will denote by $\mathbf{F}_{\Sigma^{\boldsymbol{\mathcal{A}}^{(2)}}}^{(1,2)}(\mathbf{Pth}_{\boldsymbol{\mathcal{A}}^{(2)}})$ the $\Sigma^{\boldsymbol{\mathcal{A}}}$-algebra $(\mathbf{in}^{\Sigma,(1,2)}(\mathbf{F}_{\Sigma^{\boldsymbol{\mathcal{A}}^{(2)}}}(
\mathbf{Pth}_{\boldsymbol{\mathcal{A}}^{(2)}}))$. We will call $\mathbf{F}_{\Sigma^{\boldsymbol{\mathcal{A}}^{(2)}}}^{(1,2)}(\mathbf{Pth}_{\boldsymbol{\mathcal{A}}^{(2)}})$ the $\Sigma^{\boldsymbol{\mathcal{A}}}$-reduct of the free $\Sigma^{\boldsymbol{\mathcal{A}}^{(2)}}$-completion $\mathbf{F}_{\Sigma^{\boldsymbol{\mathcal{A}}^{(2)}}}(
\mathbf{Pth}_{\boldsymbol{\mathcal{A}}^{(2)}})$.
\end{definition}

\begin{figure}
\begin{tikzpicture}
[ACliment/.style={-{To [angle'=45, length=5.75pt, width=4pt, round]},
scale=.8
}]
\node[] (X) 	at 	(0,0) 	[] 	{$X$};
\node[] (PT) 	at 	(4,0) 	[] 	{$[\mathbf{PT}_{\boldsymbol{\mathcal{A}}}]$};
\node[] (P2)	at		(4,-2)	[]	{$\mathbf{Pth}_{\boldsymbol{\mathcal{A}}^{(2)}}$};
\node[] (FP2)	at		(8,-2)	 []
{$\mathbf{F}_{\Sigma^{\boldsymbol{\mathcal{A}}^{(2)}}}(\mathbf{Pth}_{\boldsymbol{\mathcal{A}}^{(2)}})$};
\node[] (T2)	at		(4,-4)	[]	{$\mathbf{T}_{\Sigma^{\boldsymbol{\mathcal{A}}^{(2)}}}(X)$};

\draw[ACliment]  (X) 	to node [above]	{$\eta^{([1],X)}$} (PT);
\draw[ACliment, bend right=10]  (X) 	to node [below left]	 {$\mathrm{ip}^{(2,X)}$} (P2);
\draw[ACliment]  (P2) 	to node [above]	{$\eta^{(2,\mathbf{Pth}_{\boldsymbol{\mathcal{A}}^{(2)}})}$} (FP2);
\draw[ACliment]  (PT) 	to node [left]	{$\mathrm{ip}^{(2,[1])\sharp}$} (P2);
\draw[ACliment]  (P2) 	to node [left]	{$\mathrm{CH}^{(2)}$} (T2);
\draw[ACliment, bend right=10]  (T2) 	to node [below right]	 {$\mathrm{ip}^{(2,X)@}$} (FP2);
\end{tikzpicture}
\caption{$\mathrm{ip}^{(2,X)@}$ on $(2,[1])$-identity second-order paths.}
\label{FDIpDU}
\end{figure}

\begin{proposition}
\label{PDIpDU} The following equation holds
\begin{itemize}
\item[(i)] $
\mathrm{ip}^{(2,X)@}\circ\eta^{(2,1)\sharp}
=
\eta^{(2,\mathbf{Pth}_{\boldsymbol{\mathcal{A}}^{(2)}})}
\circ
\mathrm{ip}^{(2,[1])\sharp}
\circ
\mathrm{pr}^{\Theta^{[1]}};
$
\item[(ii)] $\mathrm{ip}^{(2,X)@}
\circ
\mathrm{CH}^{(2)}
\circ
\mathrm{ip}^{(2,[1])\sharp}
=
\eta^{(2,\mathbf{Pth}_{\boldsymbol{\mathcal{A}}^{(2)}})}
\circ
\mathrm{ip}^{(2,[1])\sharp}.
$
\end{itemize}
\end{proposition}
\begin{proof}
The reader is advised to consult the diagram in Figure~\ref{FDIpDU} for a better insight of the many-sorted mappings under consideration.

We first prove that the following equation holds
$$
\mathrm{ip}^{(2,X)@}\circ\eta^{(2,1)\sharp}
=
\eta^{(2,\mathbf{Pth}_{\boldsymbol{\mathcal{A}}^{(2)}})}
\circ
\mathrm{ip}^{(2,[1])\sharp}
\circ
\mathrm{pr}^{\Theta^{[1]}}
.
$$

Note that the $S$-sorted mapping $\mathrm{ip}^{(2,X)@}\circ\eta^{(2,1)\sharp}$ is a $\Sigma^{\boldsymbol{\mathcal{A}}}$-homomorphism from $\mathbf{PT}_{\boldsymbol{\mathcal{A}}}$ to $\mathbf{F}^{(1,2)}_{\Sigma^{\boldsymbol{\mathcal{A}}^{(2)}}}(\mathbf{Pth}_{\boldsymbol{\mathcal{A}}^{(2)}})$ in virtue of Definition~\ref{DDIp} and Proposition~\ref{PDUEmb}.

Moreover, the following chain of equalities holds
\allowdisplaybreaks
\begin{align*}
\mathrm{ip}^{(2,X)@}
\circ
\eta^{(2,1)\sharp}
\circ
\eta^{(1,X)}
&=
\mathrm{ip}^{(2,X)@}
\circ
\eta^{(2,X)}
\tag{1}
\\&=
\eta^{(2,\mathbf{Pth}_{\boldsymbol{\mathcal{A}}^{(2)}})}
\circ
\mathrm{ip}^{(2,X)}.
\tag{2}
\end{align*}

The first equality follows from Proposition~\ref{PDUEmb}; whilst the second equality follows from Definition~\ref{DDIp}.

On the other hand the $S$-sorted mapping $\eta^{(2,\mathbf{Pth}_{\boldsymbol{\mathcal{A}}^{(2)}})}\circ\mathrm{ip}^{(2,[1])\sharp}\circ\mathrm{pr}^{\Theta^{[1]}}$ is a $\Sigma^{\boldsymbol{\mathcal{A}}}$-homomorphism from $\mathbf{T}_{\Sigma^{\boldsymbol{\mathcal{A}}}}(X)$ to $\mathbf{F}^{(1,2)}_{\Sigma^{\boldsymbol{\mathcal{A}}^{(2)}}}(\mathbf{Pth}_{\boldsymbol{\mathcal{A}}^{(2)}})$ in virtue of Remark~\ref{RDFP}, Proposition~\ref{PDUIpCatHom} and Corollary~\ref{CPTQPr}.

Moreover, the following equality holds
$$
\eta^{(2,\mathbf{Pth}_{\boldsymbol{\mathcal{A}}^{(2)}})}
\circ
\mathrm{ip}^{(2,[1])\sharp}
\circ
\mathrm{pr}^{\Theta^{[1]}}
\circ
\eta^{(1,X)}
=
\eta^{(2,\mathbf{Pth}_{\boldsymbol{\mathcal{A}}^{(2)}})}
\circ
\mathrm{ip}^{(2,X)}.
$$
This is so, in virtue of Proposition~\ref{PDBasicEq} and Definition~\ref{DPTQEta}.

By the Universal Property of the free $\Sigma^{\boldsymbol{\mathcal{A}}}$-algebra $\mathbf{T}_{\Sigma^{\boldsymbol{\mathcal{A}}}}(X)$ we have the desired equality
\begin{align*}
\mathrm{ip}^{(2,X)@}\circ\eta^{(2,1)\sharp}
&=
\eta^{(2,\mathbf{Pth}_{\boldsymbol{\mathcal{A}}^{(2)}})}
\circ
\mathrm{ip}^{(2,[1])\sharp}
\circ
\mathrm{pr}^{\Theta^{[1]}}
.
\tag{$\star$}
\end{align*}

Therefore, the following chain of equalities holds
\begin{flushleft}
$\mathrm{ip}^{(2,X)@}
\circ
\mathrm{CH}^{(2)}
\circ
\mathrm{ip}^{(2,[1])\sharp}$
\allowdisplaybreaks
\begin{align*}
\qquad
&=
\mathrm{ip}^{(2,X)@}
\circ
\eta^{(2,1)\sharp}
\circ
\mathrm{CH}^{(1)\mathrm{m}}
\circ
\mathrm{ip}^{([1],X)@}
\tag{1}
\\&=
\eta^{(2,\mathbf{Pth}_{\boldsymbol{\mathcal{A}}^{(2)}})}
\circ
\mathrm{ip}^{(2,[1])\sharp}
\circ
\mathrm{pr}^{\Theta^{[1]}}
\circ
\mathrm{CH}^{(1)\mathrm{m}}
\circ
\mathrm{ip}^{([1],X)@}
\tag{2}
\\&=
\eta^{(2,\mathbf{Pth}_{\boldsymbol{\mathcal{A}}^{(2)}})}
\circ
\mathrm{ip}^{(2,[1])\sharp}.
\tag{3}
\end{align*}
\end{flushleft}

The first equality follows from Proposition~\ref{PDCHDUId}; the second equality follows from equality~($\star$) above; finally, the last equality follows from the fact that $\mathrm{pr}^{\Theta^{[1]}}
\circ
\mathrm{CH}^{(1)\mathrm{m}}
\circ
\mathrm{ip}^{([1],X)@}$ sends, for every sort $s\in S$, a path term class $[P]_{s}$ in $[\mathbf{PT}_{\boldsymbol{\mathcal{A}}}]_{s}$ to 
$
[\mathrm{CH}^{(1)}_{s}(\mathrm{ip}^{(1,X)@}_{s}(P))]_{s},
$ which, by Lemma~\ref{LWCong}, coincides with $[P]_{s}$, i.e., 
$$
\mathrm{pr}^{\Theta^{[1]}}
\circ
\mathrm{CH}^{(1)\mathrm{m}}
\circ
\mathrm{ip}^{([1],X)@}=\mathrm{id}^{
[\mathbf{PT}_{\boldsymbol{\mathcal{A}}}]
}.
$$

This finishes the proof.
\end{proof}

\begin{figure}
\begin{tikzpicture}
[ACliment/.style={-{To [angle'=45, length=5.75pt, width=4pt, round]},
scale=.8
}]
\node[] (X) 	at 	(0,0) 	[] 	{$X$};
\node[] (PT) 	at 	(4,0) 	[] 	{$\mathbf{T}_{\Sigma}(X)$};
\node[] (P2)	at		(4,-2)	[]	{$\mathbf{Pth}_{\boldsymbol{\mathcal{A}}^{(2)}}$};
\node[] (FP2)	at		(8,-2)	 []
{$\mathbf{F}_{\Sigma^{\boldsymbol{\mathcal{A}}^{(2)}}}(\mathbf{Pth}_{\boldsymbol{\mathcal{A}}^{(2)}})$};
\node[] (T2)	at		(4,-4)	[]	{$\mathbf{T}_{\Sigma^{\boldsymbol{\mathcal{A}}^{(2)}}}(X)$};

\draw[ACliment]  (X) 	to node [above]	{$\eta^{(0,X)}$} (PT);
\draw[ACliment, bend right=10]  (X) 	to node [below left]	 {$\mathrm{ip}^{(2,X)}$} (P2);
\draw[ACliment]  (P2) 	to node [above]	{$\eta^{(2,\mathbf{Pth}_{\boldsymbol{\mathcal{A}}^{(2)}})}$} (FP2);
\draw[ACliment]  (PT) 	to node [left]	{$\mathrm{ip}^{(2,0)\sharp}$} (P2);
\draw[ACliment]  (P2) 	to node [left]	{$\mathrm{CH}^{(2)}$} (T2);
\draw[ACliment, bend right=10]  (T2) 	to node [below right]	 {$\mathrm{ip}^{(2,X)@}$} (FP2);
\end{tikzpicture}
\caption{$\mathrm{ip}^{(2,X)@}$ on $(2,0)$-identity second-order paths.}
\label{FDIpDZ}
\end{figure}

We now explore the relation between the mapping $\mathrm{ip}^{(2,X)@}$ and the value of the second-order Curry-Howard mapping at  $(2,0)$-identity second-order paths. Before doing so, we introduce the $\Sigma$-reduct of the free $\Sigma^{\boldsymbol{\mathcal{A}}^{(2)}}$-completion of the many-sorted partial $\Sigma^{\boldsymbol{\mathcal{A}}^{(2)}}$-algebra of second-order paths

\begin{definition}\label{DDPFRedDZ} Let $\mathrm{in}^{\Sigma,(2,0)}$ be the canonical embedding of $\Sigma$ into $\Sigma^{\boldsymbol{\mathcal{A}}^{(2)}}$. Then by Proposition~\ref{PFunSig}, for the morphism $\mathbf{in}^{\Sigma,(2,0)}=(\mathrm{id}_{S},\mathrm{in}^{\Sigma,(2,0)})$ from $(S,\Sigma)$ to $(S,\Sigma^{\boldsymbol{\mathcal{A}}^{(2)}})$ and the free $\Sigma^{\boldsymbol{\mathcal{A}}^{(2)}}$-completion $\mathbf{F}_{\Sigma^{\boldsymbol{\mathcal{A}}^{(2)}}}(\mathbf{Pth}_{\boldsymbol{\mathcal{A}}^{(2)}})$, we will denote by $\mathbf{F}_{\Sigma^{\boldsymbol{\mathcal{A}}^{(2)}}}^{(0,2)}(\mathbf{Pth}_{\boldsymbol{\mathcal{A}}^{(2)}})$ the $\Sigma$-algebra $(\mathbf{in}^{\Sigma,(0,2)}(\mathbf{F}_{\Sigma^{\boldsymbol{\mathcal{A}}^{(2)}}}(
\mathbf{Pth}_{\boldsymbol{\mathcal{A}}^{(2)}}))$. We will call $\mathbf{F}_{\Sigma^{\boldsymbol{\mathcal{A}}^{(2)}}}^{(0,2)}(\mathbf{Pth}_{\boldsymbol{\mathcal{A}}^{(2)}})$ the $\Sigma$-reduct of the free $\Sigma^{\boldsymbol{\mathcal{A}}^{(2)}}$-completion $\mathbf{F}_{\Sigma^{\boldsymbol{\mathcal{A}}^{(2)}}}(
\mathbf{Pth}_{\boldsymbol{\mathcal{A}}^{(2)}})$.
\end{definition}

\begin{proposition}
\label{PDIpDZ} The following equations holds
\begin{itemize}
\item[(i)] $
\mathrm{ip}^{(2,X)@}\circ\eta^{(2,0)\sharp}
=
\eta^{(2,\mathbf{Pth}_{\boldsymbol{\mathcal{A}}^{(2)}})}
\circ
\mathrm{ip}^{(2,0)\sharp};
$
\item[(ii)] $\mathrm{ip}^{(2,X)@}
\circ
\mathrm{CH}^{(2)}
\circ
\mathrm{ip}^{(2,0)\sharp}
=
\eta^{(2,\mathbf{Pth}_{\boldsymbol{\mathcal{A}}^{(2)}})}
\circ
\mathrm{ip}^{(2,0)\sharp}.
$
\end{itemize}
\end{proposition}
\begin{proof}
The reader is advised to consult the diagram in Figure~\ref{FDIpDZ} for a better insight of the many-sorted mappings under consideration.

We first prove that the following equation holds
$$
\mathrm{ip}^{(2,X)@}\circ\eta^{(2,0)\sharp}
=
\eta^{(2,\mathbf{Pth}_{\boldsymbol{\mathcal{A}}^{(2)}})}
\circ
\mathrm{ip}^{(2,0)\sharp}
.
$$

Note that the $S$-sorted mapping $\mathrm{ip}^{(2,X)@}\circ\eta^{(2,0)\sharp}$ is a $\Sigma$-homomorphism from $\mathbf{T}_{\Sigma}(X)$ to $\mathbf{F}^{(0,2)}_{\Sigma^{\boldsymbol{\mathcal{A}}^{(2)}}}(\mathbf{Pth}_{\boldsymbol{\mathcal{A}}^{(2)}})$ in virtue of Definition~\ref{DDIp} and Proposition~\ref{PDZEmb}.

Moreover, the following chain of equalities holds
\begin{align*}
\mathrm{ip}^{(2,X)@}
\circ
\eta^{(2,0)\sharp}
\circ
\eta^{(0,X)}
&=
\mathrm{ip}^{(2,X)@}
\circ
\eta^{(2,X)}
\tag{1}
\\&=
\eta^{(2,\mathbf{Pth}_{\boldsymbol{\mathcal{A}}^{(2)}})}
\circ
\mathrm{ip}^{(2,X)}.
\tag{2}
\end{align*}

The first equality follows from Proposition~\ref{PDZEmb}; whilst the second equality follows from Definition~\ref{DDIp}.

On the other hand the $S$-sorted mapping $\eta^{(2,\mathbf{Pth}_{\boldsymbol{\mathcal{A}}^{(2)}})}\circ\mathrm{ip}^{(2,0)\sharp}$ is a $\Sigma$-homomorphism from $\mathbf{T}_{\Sigma}(X)$ to $\mathbf{F}^{(0,2)}_{\Sigma^{\boldsymbol{\mathcal{A}}^{(2)}}}(\mathbf{Pth}_{\boldsymbol{\mathcal{A}}^{(2)}})$ in virtue of Remark~\ref{RDFP} and Proposition~\ref{PDZHomIp}.

Moreover, the following equality holds
$$
\eta^{(2,\mathbf{Pth}_{\boldsymbol{\mathcal{A}}^{(2)}})}
\circ
\mathrm{ip}^{(2,0)\sharp}
\circ
\eta^{(0,X)}
=
\eta^{(2,\mathbf{Pth}_{\boldsymbol{\mathcal{A}}^{(2)}})}
\circ
\mathrm{ip}^{(2,X)}.
$$
This is so, in virtue of Proposition~\ref{PDBasicEqZ}.

By the Universal Property of the free $\Sigma$-algebra $\mathbf{T}_{\Sigma}(X)$ we have the desired equality
\begin{align*}
\mathrm{ip}^{(2,X)@}\circ\eta^{(2,0)\sharp}
&=
\eta^{(2,\mathbf{Pth}_{\boldsymbol{\mathcal{A}}^{(2)}})}
\circ
\mathrm{ip}^{(2,0)\sharp}
.
\tag{$\star$}
\end{align*}

Therefore, the following chain of equalities holds
\allowdisplaybreaks
\begin{align*}
\mathrm{ip}^{(2,X)@}
\circ
\mathrm{CH}^{(2)}
\circ
\mathrm{ip}^{(2,0)\sharp}
&=
\mathrm{ip}^{(2,X)@}
\circ
\eta^{(2,0)\sharp}
\tag{1}
\\&=
\eta^{(2,\mathbf{Pth}_{\boldsymbol{\mathcal{A}}^{(2)}})}
\circ
\mathrm{ip}^{(2,0)\sharp}.
\tag{2}
\end{align*}

The first equality follows from Proposition~\ref{PDCHDUId}; finally, the last equality follows from equality~($\star$) above.

This finishes the proof.
\end{proof}

\begin{proposition}
\label{PDIpEch} The following equations hold
\begin{itemize}
\item[(i)] $
\mathrm{ip}^{(2,X)@}\circ\eta^{(2,\mathcal{A})}
=
\mathrm{ech}^{(2,\mathcal{A})};
$
\item[(ii)] $
\mathrm{ip}^{(2,X)@}\circ\eta^{(2,\mathcal{A}^{(2)})}
=
\mathrm{ech}^{(2,\mathcal{A}^{(2)})}.
$
\end{itemize}
\end{proposition}
\begin{proof}
Let $s$ be a sort in $S$ and let $\mathfrak{p}$ be a rewrite rule in $\mathcal{A}_{s}$.

The following chain of equalities holds
\allowdisplaybreaks
\begin{align*}
\mathrm{ip}^{(2,X)@}_{s}\left(\eta^{(2,\mathcal{A})}_{s}\left(
\mathfrak{p}
\right)\right)
&=
\mathrm{ip}^{(2,X)@}_{s}\left(
\mathfrak{p}^{\mathbf{T}_{\Sigma^{\boldsymbol{\mathcal{A}}^{(2)}}}(X)}
\right)
\tag{1}
\\&=
\mathfrak{p}^{
\mathbf{F}_{\Sigma^{\boldsymbol{\mathcal{A}}^{(2)}}}(\mathbf{Pth}_{\boldsymbol{\mathcal{A}}^{(2)}})
}
\tag{2}
\\&=
\mathfrak{p}^{
\mathbf{Pth}_{\boldsymbol{\mathcal{A}}^{(2)}}
}
\tag{3}
\\&=
\mathrm{ech}^{(2,\mathcal{A})}_{s}\left(
\mathfrak{p}
\right).
\tag{4}
\end{align*}

In the just stated chain of equalities, the first equality unravels the value of the mapping $\eta^{(2,\mathcal{A})}$ at $\mathfrak{p}$, according to Definition~\ref{DDPth}; the second equality follows from the fact that $\mathrm{ip}^{(2,X)@}$ is a $\Sigma^{\boldsymbol{\mathcal{A}}^{(2)}}$-homomorphism, according to Definition~\ref{DDIp}; the third equality follows from the fact that there is an interpretation for the constant operation symbol $\mathfrak{p}$ in the many-sorted partial $\Sigma^{\boldsymbol{\mathcal{A}}^{(2)}}$-algebra $\mathbf{Pth}_{\boldsymbol{\mathcal{A}}^{(2)}}$, according to Proposition~\ref{PDPthCatAlg}, thus the interpretation of the constant operation symbol $\mathfrak{p}$ in the many-sorted  $\Sigma^{\boldsymbol{\mathcal{A}}^{(2)}}$-algebra $\mathbf{F}_{\Sigma^{\boldsymbol{\mathcal{A}}^{(2)}}}(\mathbf{Pth}_{\boldsymbol{\mathcal{A}}^{(2)}})$ becomes that of $\mathbf{Pth}_{\boldsymbol{\mathcal{A}}^{(2)}}$; finally, the last equality recovers the description of the constant operation symbol $\mathfrak{p}$ in the many-sorted partial $\Sigma^{\boldsymbol{\mathcal{A}}^{(2)}}$-algebra $\mathbf{Pth}_{\boldsymbol{\mathcal{A}}^{(2)}}$, according to Proposition~\ref{PDPthCatAlg}.

Now, let $s$ be a sort in $S$ and let $\mathfrak{p}^{(2)}$ be a second-order rewrite rule in $\mathcal{A}^{(2)}_{s}$.

The following chain of equalities holds
\allowdisplaybreaks
\begin{align*}
\mathrm{ip}^{(2,X)@}_{s}\left(\eta^{(2,\mathcal{A}^{(2)})}_{s}\left(
\mathfrak{p}^{(2)}
\right)\right)
&=
\mathrm{ip}^{(2,X)@}_{s}\left(
\mathfrak{p}^{(2)\mathbf{T}_{\Sigma^{\boldsymbol{\mathcal{A}}^{(2)}}}(X)}
\right)
\tag{1}
\\&=
\mathfrak{p}^{(2)
\mathbf{F}_{\Sigma^{\boldsymbol{\mathcal{A}}^{(2)}}}(\mathbf{Pth}_{\boldsymbol{\mathcal{A}}^{(2)}})
}
\tag{2}
\\&=
\mathfrak{p}^{(2)
\mathbf{Pth}_{\boldsymbol{\mathcal{A}}^{(2)}}
}
\tag{3}
\\&=
\mathrm{ech}^{(2,\mathcal{A}^{(2)})}_{s}\left(
\mathfrak{p}^{(2)}
\right).
\tag{4}
\end{align*}

In the just stated chain of equalities, the first equality unravels the value of the mapping $\eta^{(2,\mathcal{A}^{(2)})}$ at $\mathfrak{p}^{(2)}$, according to Definition~\ref{DDEta}; the second equality follows from the fact that $\mathrm{ip}^{(2,X)@}$ is a $\Sigma^{\boldsymbol{\mathcal{A}}^{(2)}}$-homomorphism, according to Definition~\ref{DDIp}; the third equality follows from the fact that there is an interpretation for the constant operation symbol $\mathfrak{p}^{(2)}$ in the many-sorted partial $\Sigma^{\boldsymbol{\mathcal{A}}^{(2)}}$-algebra $\mathbf{Pth}_{\boldsymbol{\mathcal{A}}^{(2)}}$, according to Proposition~\ref{PDPthDCatAlg}, thus the interpretation of the constant operation symbol $\mathfrak{p}^{(2)}$ in the many-sorted  $\Sigma^{\boldsymbol{\mathcal{A}}^{(2)}}$-algebra $\mathbf{F}_{\Sigma^{\boldsymbol{\mathcal{A}}^{(2)}}}(\mathbf{Pth}_{\boldsymbol{\mathcal{A}}^{(2)}})$ becomes that of $\mathbf{Pth}_{\boldsymbol{\mathcal{A}}^{(2)}}$; finally, the last equality recovers the description of the constant operation symbol $\mathfrak{p}^{(2)}$ in the many-sorted partial $\Sigma^{\boldsymbol{\mathcal{A}}^{(2)}}$-algebra $\mathbf{Pth}_{\boldsymbol{\mathcal{A}}^{(2)}}$, according to Proposition~\ref{PDPthDCatAlg}.

This completes the proof.
\end{proof}

The next proposition is fundamental for the rest of this work. It states that the $\mathrm{ip}^{(2,X)@}$ mapping acting on the value of a  second-order Curry-Howard mapping of a second-order path is always another second-order path, not necessarily equal to the original input, but in the same class for the kernel of the second-order Curry-Howard mapping. This proposition takes into account the embedding of $\mathbf{Pth}_{\boldsymbol{\mathcal{A}}^{(2)}}$ as a weak subalgebra of $\mathbf{F}_{\Sigma^{\boldsymbol{\mathcal{A}}^{(2)}}}(\mathbf{Pth}_{\boldsymbol{\mathcal{A}}^{(2)}})$, as stated in Definition~\ref{DDFP}. 

\begin{restatable}{proposition}{PDIpDCH}
\label{PDIpDCH} The mapping 
$$
\mathrm{ip}^{(2,X)@}
\circ
\mathrm{CH}^{(2)}
\colon
\mathrm{Pth}_{\boldsymbol{\mathcal{A}}^{(2)}}
\mor
\mathrm{F}_{\Sigma^{\boldsymbol{\mathcal{A}}^{(2)}}}(
\mathbf{Pth}_{\boldsymbol{\mathcal{A}}^{(2)}}
)
$$
sends, for every sort $s\in S$, a second-order path $\mathfrak{P}^{(2)}$ in $\mathrm{Pth}_{\boldsymbol{\mathcal{A}}^{(2)},s}$ to a second-order path in the equivalence class $[\mathfrak{P}^{(2)}]^{}_{s}$.
\end{restatable}
\begin{proof}

We want to show that, for every sort $s\in S$ and every second-order path $\mathfrak{P}^{(2)}$ in $\mathrm{Pth}_{\boldsymbol{\mathcal{A}}^{(2)},s}$, the element $\mathrm{ip}^{(2,X)@}_{s}(\mathrm{CH}^{(2)}_{s}(\mathfrak{P}^{(2)}))$ of $\mathrm{F}_{\Sigma^{\boldsymbol{\mathcal{A}}^{(2)}}}(\mathbf{Pth}_{\boldsymbol{\mathcal{A}}^{(2)}})_{s}$ satisfies:
\begin{itemize}
\item[(i)] $\mathrm{ip}^{(2,X)@}_{s}(\mathrm{CH}^{(2)}_{s}(\mathfrak{P}^{(2)}))$ is a second-order path in $\mathrm{Pth}_{\boldsymbol{\mathcal{A}}^{(2)},s}$; and
\item[(ii)] $\mathrm{ip}^{(2,X)@}_{s}(\mathrm{CH}^{(2)}_{s}(\mathfrak{P}^{(2)}))$ belongs to $[\mathfrak{P}^{(2)}]^{}_{s}$.
\end{itemize}

In case of $\mathfrak{P}^{(2)}$ being a $(2,[1])$-identity second-order path, these results follow from Proposition~\ref{PDIpDU}. Note that, for every sort $s\in S$, if $\mathfrak{P}^{(2)}$ is a $(2,[1])$-identity second-order path, then there exists a path term class $[P]_{s}$ in $[\mathbf{PT}_{\boldsymbol{\mathcal{A}}}]_{s}$,
for which $\mathfrak{P}^{(2)}=\mathrm{ip}^{(2,[1])\sharp}_{s}([P]_{s})$.

The following chain of equalities holds
\allowdisplaybreaks
\begin{align*}
\mathrm{ip}^{(2,X)@}_{s}\left(\mathrm{CH}^{(2)}_{s}\left(
\mathfrak{P}^{(2)}
\right)\right)
&=
\mathrm{ip}^{(2,X)@}_{s}\left(\mathrm{CH}^{(2)}_{s}\left(
\mathrm{ip}^{(2,[1])\sharp}_{s}\left(
[P]_{s}
\right)
\right)\right)
\tag{1}
\\&=
\eta^{(2,\mathbf{Pth}_{\boldsymbol{\mathcal{A}}^{(2)}})}_{s}\left(
\mathrm{ip}^{(2,[1])\sharp}_{s}\left([P]_{s}
\right)\right)
\tag{2}
\\&=
\eta^{(2,\mathbf{Pth}_{\boldsymbol{\mathcal{A}}^{(2)}})}_{s}\left(
\mathfrak{P}^{(2)}
\right)
.
\tag{3}
\end{align*}

The first equality unravels the description of $\mathfrak{P}^{(2)}$ as a $(2,[1])$-identity second-order path; the second equality follows from Proposition~\ref{PDIpDU}; finally, the last equality recovers the description of $\mathfrak{P}^{(2)}$ as a $(2,[1])$-identity second-order path.

We identify $\eta^{(2,\mathbf{Pth}_{\boldsymbol{\mathcal{A}}^{(2)}})}_{s}(
\mathfrak{P}^{(2)}
)$ with $\mathfrak{P}^{(2)}$. Note that  
\begin{itemize}
\item[(i)] $\mathrm{ip}^{(2,X)@}_{s}(\mathrm{CH}^{(2)}_{s}(
\mathfrak{P}^{(2)}
))$  is a second-order path in $\mathrm{Pth}_{\boldsymbol{\mathcal{A}}^{(2)},s}$, and
\item[(ii)] $\mathrm{ip}^{(2,X)@}_{s}(\mathrm{CH}^{(2)}_{s}(
\mathfrak{P}^{(2)}
))$  belongs to $[\mathfrak{P}^{(2)}]^{}_{s}$.
\end{itemize}

We can, therefore, assume without loss of generality that $\mathfrak{P}^{(2)}$ is a second-order path of length at least one.

We prove the general case by Artinian induction on 
$(
\coprod\mathrm{Pth}_{\boldsymbol{\mathcal{A}}^{(2)}},
\leq_{\mathbf{Pth}_{\boldsymbol{\mathcal{A}}^{(2)}}}
)
$.

\textsf{Base step of the Artinian induction.}

Let $(\mathfrak{P}^{(2)},s)$ be a minimal element in $(\coprod\mathrm{Pth}_{\boldsymbol{\mathcal{A}}^{(2)}},\leq_{\mathbf{Pth}_{\boldsymbol{\mathcal{A}}^{(2)}}})
$. Then by Proposition~\ref{PDMinimal}, and taking into account that we are assuming that $\mathfrak{P}^{(2)}$ is a non-$(2,[1])$-identity second-order path, we have that $\mathfrak{P}^{(2)}$ is a  second-order echelon associated to a second-order rewrite rule $\mathfrak{p}^{(2)}\in\mathcal{A}^{(2)}_{s}$. 

According to Definition~\ref{DDCH}, the value of the second-order Curry-Howard mapping at $\mathfrak{P}^{(2)}$ is given by 
$$
\mathrm{CH}^{(2)}_{s}\left(
\mathfrak{P}^{(2)}
\right)
=
\mathfrak{p}^{(2)
\mathbf{T}_{\Sigma^{\boldsymbol{\mathcal{A}}^{(2)}}}(X)
}.
$$

Thus we have that 
\begin{align*}
\mathrm{ip}^{(2,X)@}_{s}\left(
\mathrm{CH}^{(2)}_{s}\left(
\mathfrak{P}^{(2)}
\right)\right)
&=
\mathrm{ip}^{(2,X)@}_{s}\left(
\mathfrak{p}^{(2)\mathbf{T}_{\Sigma^{\boldsymbol{\mathcal{A}}^{(2)}}}(X)}
\right)
\tag{1}
\\&=
\mathfrak{p}^{(2)
\mathbf{F}_{\Sigma^{\boldsymbol{\mathcal{A}}^{(2)}}}\left(
\mathbf{Pth}_{\boldsymbol{\mathcal{A}}^{(2)}}
\right)
}
\tag{2}
\\&=
\mathfrak{p}^{(2)
\mathbf{Pth}_{\boldsymbol{\mathcal{A}}^{(2)}}
}.
\tag{3}
\end{align*}

In the just stated chain of equalities, the first equality unravels the second-order Curry-Howard mapping according to Definition~\ref{DDCH}; the second equality holds because, according to Definition~\ref{DDIp}, $\mathrm{ip}^{(2,X)@}$ is a $\Sigma^{\boldsymbol{\mathcal{A}}^{(2)}}$-homomorphism; and the third equality holds since, by Proposition~\ref{PDPthAlg},  we have an interpretation of the constant symbol $\mathfrak{p}^{(2)}$ in the many-sorted partial $\Sigma^{\boldsymbol{\mathcal{A}}^{(2)}}$-algebra $\mathbf{Pth}_{\boldsymbol{\mathcal{A}}^{(2)}}$.

Let us recall that $\mathfrak{p}^{(2)
\mathbf{Pth}_{\boldsymbol{\mathcal{A}}^{(2)}}
}$ is the second-order echelon
$\mathrm{ech}^{(2,\mathcal{A}^{(2)})}_{s}(\mathfrak{p}^{(2)})$, i.e., the original second-order path $\mathfrak{P}^{(2)}$. 

Therefore, 
\begin{itemize}
\item[(i)] $\mathrm{ip}^{(2,X)@}_{s}(\mathrm{CH}^{(2)}_{s}(
\mathfrak{P}^{(2)}
))$  is a second-order path in $\mathrm{Pth}_{\boldsymbol{\mathcal{A}}^{(2)},s}$, and
\item[(ii)] $\mathrm{ip}^{(2,X)@}_{s}(\mathrm{CH}^{(2)}_{s}(
\mathfrak{P}^{(2)}
))$  belongs to $[\mathfrak{P}^{(2)}]^{}_{s}$.
\end{itemize}

This concludes the base case.

\textsf{Inductive step of the Artinian induction.}

Let $(\mathfrak{P}^{(2)},s)$ be a non-minimal element in $(\coprod\mathrm{Pth}_{\boldsymbol{\mathcal{A}}^{(2)}}, \leq_{\mathbf{Pth}_{\boldsymbol{\mathcal{A}}^{(2)}}})$. Let us suppose that, for every sort $t\in S$ and every second-order path $\mathfrak{Q}^{(2)}$ in $\mathrm{Pth}_{\boldsymbol{\mathcal{A}}^{(2)},t}$, if $(\mathfrak{Q}^{(2)},t)<_{\mathbf{Pth}_{\boldsymbol{\mathcal{A}}^{(2)}}} (\mathfrak{P}^{(2)},s)$, then the statement holds for $\mathfrak{Q}^{(2)}$, i.e., 
\begin{itemize}
\item[(i)] $\mathrm{ip}^{(2,X)@}_{t}(\mathrm{CH}^{(2)}_{t}(\mathfrak{Q}^{(2)}))$ is a second-order path in $\mathrm{Pth}_{\boldsymbol{\mathcal{A}}^{(2)},t}$; and
\item[(ii)] $\mathrm{ip}^{(2,X)@}_{t}(\mathrm{CH}^{(2)}_{t}(\mathfrak{Q}^{(2)}))$ belongs to $[\mathfrak{Q}^{(2)}]^{}_{t}$.
\end{itemize}

Since $(\mathfrak{P}^{(2)},s)$ is a non-minimal element in $(\coprod\mathrm{Pth}_{\boldsymbol{\mathcal{A}}^{(2)}},\leq_{\mathbf{Pth}_{\boldsymbol{\mathcal{A}}^{(2)}}})$ and taking into account that $\mathfrak{P}^{(2)}$ is a non-$(2,[1])$-identity second-order path, we have, by Lemma~\ref{LDOrdI}, that $\mathfrak{P}^{(2)}$ is either (1) a second-order path of length strictly greater than one containing at least one  second-order echelon or (2) an echelonless second-order path.

If (1), i.e., if $\mathfrak{P}^{(2)}$ is a second-order path of length strictly greater than one containing at least one  second-order echelon. Then, we let $i\in\bb{\mathfrak{P}^{(2)}}$ be the first index for which the one-step subpath $\mathfrak{P}^{(2),i,i}$ is a  second-order echelon. We distinguish two cases accordingly, either (1.1) $i=0$ or~(1.2) $i>0$.

If~(1.1), i.e., if $\mathfrak{P}^{(2)}$  is a second-order path of length strictly greater than one having its first  second-order echelon on its first step then, according to Definition~\ref{DDCH}, the value of the second-order Curry-Howard mapping at $\mathfrak{P}^{(2)}$ is given by
$$
\mathrm{CH}^{(2)}_{s}\left(
\mathfrak{P}^{(2)}
\right)
=
\mathrm{CH}^{(2)}_{s}\left(
\mathfrak{P}^{(2),1,\bb{\mathfrak{P}^{(2)}}-1}
\right)
\circ_{s}^{1\mathbf{T}_{\Sigma^{\boldsymbol{\mathcal{A}}^{(2)}}}(X)}
\mathrm{CH}^{(2)}_{s}\left(
\mathfrak{P}^{(2),0,0}
\right).
$$

Therefore, the  following chain of equalities holds
\begin{flushleft}
$
\mathrm{ip}^{(2,X)@}_{s}\left(
\mathrm{CH}^{(2)}_{s}\left(
\mathfrak{P}^{(2)}
\right)\right)$
\allowdisplaybreaks
\begin{align*}
\quad
&=
\mathrm{ip}^{(2,X)@}_{s}\left(
\mathrm{CH}^{(2)}_{s}\left(
\mathfrak{P}^{(2),1,\bb{\mathfrak{P}^{(2)}}-1}
\right)
\circ_{s}^{1\mathbf{T}_{\Sigma^{\boldsymbol{\mathcal{A}}^{(2)}}}(X)}
\mathrm{CH}^{(2)}_{s}\left(
\mathfrak{P}^{(2),0,0}
\right)
\right)
\tag{1}
\\&=
\mathrm{ip}^{(2,X)@}_{s}\left(
\mathrm{CH}^{(2)}_{s}\left(
\mathfrak{P}^{(2),1,\bb{\mathfrak{P}^{(2)}}-1}
\right)\right)
\circ_{s}^{1
\mathbf{F}_{\Sigma^{\boldsymbol{\mathcal{A}}^{(2)}}}\left(
\mathbf{Pth}_{\boldsymbol{\mathcal{A}}^{(2)}}
\right)}
\\&
\qquad\qquad\qquad\qquad\qquad\qquad\qquad\qquad\qquad\qquad
\mathrm{ip}^{(2,X)@}_{s}\left(
\mathrm{CH}^{(2)}_{s}\left(
\mathfrak{P}^{(2),0,0}
\right)\right)
\tag{2}
\\&=
\mathrm{ip}^{(2,X)@}_{s}\left(
\mathrm{CH}^{(2)}_{s}\left(
\mathfrak{P}^{(2),1,\bb{\mathfrak{P}^{(2)}}-1}
\right)\right)
\circ_{s}^{1
\mathbf{Pth}_{\boldsymbol{\mathcal{A}}^{(2)}}
}
\mathrm{ip}^{(2,X)@}_{s}\left(
\mathrm{CH}^{(2)}_{s}\left(
\mathfrak{P}^{(2),0,0}
\right)\right).
\tag{3}
\end{align*}
\end{flushleft}

In the just stated chain of equalities, the first equality recovers the value of the second-order Curry-Howard mapping at $\mathfrak{P}^{(2)}$; the second equality holds since, according to Definition~\ref{DDIp}, $\mathrm{ip}^{(2,X)@}$ is a many-sorted $\Sigma^{\boldsymbol{\mathcal{A}}^{(2)}}$-homomorphism; finally, the last equality follows taking into account that $(\mathfrak{P}^{(2),1,\bb{\mathfrak{P}^{(2)}}-1},s)$ and $(\mathfrak{P}^{(2),0,0},s)$ are strictly smaller than $(\mathfrak{P}^{(2)},s)$ with respect to $\prec_{\mathbf{Pth}_{\boldsymbol{\mathcal{A}}^{(2)}}}$. Then, by the inductive hypothesis, we have that $\mathrm{ip}^{(2,X)@}_{s}(\mathrm{CH}^{(2)}_{s}(\mathfrak{P}^{(2),1,\bb{\mathfrak{P}^{(2)}}-1}))$ and $\mathrm{ip}^{(2,X)@}_{s}(\mathrm{CH}^{(2)}_{s}(\mathfrak{P}^{(2),0,0}))$ are second-order paths in $[\mathfrak{P}^{(2),1,\bb{\mathfrak{P}^{(2)}}-1}]^{}_{s}$ and $[\mathfrak{P}^{(2),0,0}]^{}_{s}$ respectively. Hence the interpretation of the $1$-composition operation symbol $\circ^{1}_{s}$ in $\mathbf{F}_{\Sigma^{\boldsymbol{\mathcal{A}}^{(2)}}}(\mathbf{Pth}_{\boldsymbol{\mathcal{A}}^{(2)}})$ becomes that of $\mathbf{Pth}_{\boldsymbol{\mathcal{A}}^{(2)}}$. Note that the $1$-composition is defined because second-order paths in the same class have the same $([1],2)$-source and $([1],2)$-target according to Lemma~\ref{LDCH}. 

Hence $\mathrm{ip}^{(2,X)@}_{s}(\mathrm{CH}^{(2)}_{s}(\mathfrak{P}^{(2)}))$ is a second-order path.

Let us also note that, since $\mathrm{ip}^{(2,X)@}_{s}(\mathrm{CH}^{(2)}_{s}(\mathfrak{P}^{(2),0,0}))$ is a second-order path in $[\mathfrak{P}^{(2),0,0}]^{}_{s}$, we have, by Lemma~\ref{LDCHDEch}, that $\mathrm{ip}^{(2,X)@}_{s}(\mathrm{CH}^{(2)}_{s}(\mathfrak{P}^{(2),0,0}))$ is a  second-order echelon. Moreover, since $\mathrm{ip}^{(2,X)@}_{s}(\mathrm{CH}^{(2)}_{s}(\mathfrak{P}^{(2),1,\bb{\mathfrak{P}^{(2)}}-1}))$ is a second-order path in $[\mathfrak{P}^{(2),1,\bb{\mathfrak{P}^{(2)}}-1}]^{}_{s}$, we have, by Lemma~\ref{LDCH}, that $\mathrm{ip}^{(2,X)@}_{s}(\mathrm{CH}^{(2)}_{s}(\mathfrak{P}^{(2),1,\bb{\mathfrak{P}^{(2)}}-1}))$ is a second-order path of length equal to $\bb{\mathfrak{P}^{(2),1,\bb{\mathfrak{P}^{(2)}}-1}}$. 

All in all, we can affirm that $\mathrm{ip}^{(2,X)@}_{s}(\mathrm{CH}^{(2)}_{s}(\mathfrak{P}^{(2)}))$ is a second-order path of length strictly greater than one containing a  second-order echelon on its first step. Hence, according to Definition~\ref{DDCH}, the value of the second-order Curry-Howard mapping at $\mathrm{ip}^{(2,X)@}_{s}(\mathrm{CH}^{(2)}_{s}(\mathfrak{P}^{(2)}))$ is given by
\begin{flushleft}
$
\mathrm{CH}^{(2)}_{s}\left(
\mathrm{ip}^{(2,X)@}_{s}\left(
\mathrm{CH}^{(2)}_{s}\left(
\mathfrak{P}^{(2)}
\right)\right)\right)
$
\begin{align*}
\quad&=
\mathrm{CH}^{(2)}_{s}\left(
\left(
\mathrm{ip}^{(2,X)@}_{s}\left(
\mathrm{CH}^{(2)}_{s}\left(
\mathfrak{P}^{(2)}
\right)\right)
\right)^{1,\bb{
\mathrm{ip}^{(2,X)@}_{s}(
\mathrm{CH}^{(2)}_{s}(
\mathfrak{P}^{(2)}
))
}-1}
\right)
\circ^{1\mathbf{T}_{\Sigma^{\boldsymbol{\mathcal{A}}^{(2)}}}(X)}_{s}
\\&\qquad\qquad\qquad\qquad\qquad\qquad\qquad
\mathrm{CH}^{(2)}_{s}\left(
\left(
\mathrm{ip}^{(2,X)@}_{s}\left(
\mathrm{CH}^{(2)}_{s}\left(
\mathfrak{P}^{(2)}
\right)\right)
\right)^{0,0}\right)
\tag{1}
\\&=
\mathrm{CH}^{(2)}_{s}\left(
\mathrm{ip}^{(2,X)@}_{s}\left(
\mathrm{CH}^{(2)}_{s}\left(
\mathfrak{P}^{(2),1,\bb{\mathfrak{P}^{(2)}}-1}
\right)\right)\right)
\circ_{s}^{1
\mathbf{T}_{\Sigma^{\boldsymbol{\mathcal{A}}^{(2)}}}(X)
}
\\&\qquad\qquad\qquad\qquad\qquad\qquad\qquad\qquad
\mathrm{CH}^{(2)}_{s}\left(
\mathrm{ip}^{(2,X)@}_{s}\left(
\mathrm{CH}^{(2)}_{s}\left(
\mathfrak{P}^{(2),0,0}
\right)\right)\right)
\tag{2}
\\&=
\mathrm{CH}^{(2)}_{s}\left(
\mathfrak{P}^{(2),1,\bb{\mathfrak{P}^{(2)}}-1}
\right)
\circ_{s}^{1
\mathbf{T}_{\Sigma^{\boldsymbol{\mathcal{A}}^{(2)}}}(X)
}
\mathrm{CH}^{(2)}_{s}\left(
\mathfrak{P}^{(2),0,0}
\right)
\tag{3}
\\&=
\mathrm{CH}^{(2)}_{s}\left(
\mathfrak{P}^{(2)}
\right).
\tag{4}
\end{align*}
\end{flushleft}

The first equality follows from the fact that $\mathrm{ip}^{(2,X)@}_{s}(\mathrm{CH}^{(2)}_{s}(\mathfrak{P}^{(2)}))$ is a second-order path of length strictly greater than one containing a  second-order echelon on its first step and Definition~\ref{DDCH}; the second equality follows from the fact that, from the previous discussion, $\mathrm{ip}^{(2,X)@}_{s}(\mathrm{CH}^{(2)}_{s}(\mathfrak{P}^{(2),0,0}))$ is a  second-order echelon. Hence, from the fact that $\mathrm{ip}^{(2,X)@}_{s}(\mathrm{CH}^{(2)}_{s}(\mathfrak{P}^{(2)}))$ can be written as the $1$-composition
$$
\mathrm{ip}^{(2,X)@}_{s}\left(
\mathrm{CH}^{(2)}_{s}\left(
\mathfrak{P}^{(2),1,\bb{\mathfrak{P}^{(2)}}-1}
\right)\right)
\circ_{s}^{1\mathbf{Pth}_{\boldsymbol{\mathcal{A}}^{(2)}}}
\mathrm{ip}^{(2,X)@}_{s}\left(
\mathrm{CH}^{(2)}_{s}\left(
\mathfrak{P}^{(2),0,0}
\right)\right),
$$
we infer that 
\begin{align*}
&
\mathrm{ip}^{(2,X)@}_{s}\left(
\mathrm{CH}^{(2)}_{s}\left(
\mathfrak{P}^{(2)}
\right)\right)^{0,0}
=
\mathrm{ip}^{(2,X)@}_{s}\left(
\mathrm{CH}^{(2)}_{s}\left(
\mathfrak{P}^{(2),0,0}
\right)\right),  \mbox{ and} 
\\&
\mathrm{ip}^{(2,X)@}_{s}\left(
\mathrm{CH}^{(2)}_{s}\left(
\mathfrak{P}^{(2)}
\right)\right)^{1,\bb{\mathrm{ip}^{(2,X)@}_{s}(
\mathrm{CH}^{(2)}_{s}(
\mathfrak{P}^{(2)}
))}-1}
=
\mathrm{ip}^{(2,X)@}_{s}\left(
\mathrm{CH}^{(2)}_{s}\left(
\mathfrak{P}^{(2),1,\bb{\mathfrak{P}^{(2)}}-1}
\right)\right);
\end{align*}
the third equality follows from the fact that $\mathrm{ip}^{(2,X)@}_{s}(\mathrm{CH}^{(2)}_{s}(\mathfrak{P}^{(2),1,\bb{\mathfrak{P}^{(2)}}-1}))$ and $\mathrm{ip}^{(2,X)@}_{s}(\mathrm{CH}^{(2)}_{s}(\mathfrak{P}^{(2),0,0}))$ are second-order paths in $[\mathfrak{P}^{(2),1,\bb{\mathfrak{P}^{(2)}}-1}]^{}_{s}$
and $[\mathfrak{P}^{(2),0,0}]^{}_{s}$ respectively; finally, the last equality recovers the value of the second-order Curry-Howard mapping at $\mathfrak{P}^{(2)}$.

Therefore,
\begin{itemize}
\item[(i)] $\mathrm{ip}^{(2,X)@}_{s}(\mathrm{CH}^{(2)}_{s}(\mathfrak{P}^{(2)}))$ is a second-order path in $\mathrm{Pth}_{\boldsymbol{\mathcal{A}}^{(2)},s}$; and
\item[(ii)] $\mathrm{ip}^{(2,X)@}_{s}(\mathrm{CH}^{(2)}_{s}(\mathfrak{P}^{(2)}))$ belongs to $[\mathfrak{P}^{(2)}]^{}_{s}$.
\end{itemize}

The case $i=0$ follows.

If~(1.2), i.e., if $\mathfrak{P}^{(2)}$  is a second-order path of length strictly greater than one having its first  second-order echelon at position $i\in\bb{\mathfrak{P}^{(2)}}$, with $i>0$, then, according to Definition~\ref{DDCH}, the value of the second-order Curry-Howard mapping at $\mathfrak{P}^{(2)}$ is given by
$$
\mathrm{CH}^{(2)}_{s}\left(
\mathfrak{P}^{(2)}
\right)
=
\mathrm{CH}^{(2)}_{s}\left(
\mathfrak{P}^{(2),i,\bb{\mathfrak{P}^{(2)}}-1}
\right)
\circ_{s}^{1\mathbf{T}_{\Sigma^{\boldsymbol{\mathcal{A}}^{(2)}}}(X)}
\mathrm{CH}^{(2)}_{s}\left(
\mathfrak{P}^{(2),0,i-1}
\right).
$$

This case follows from a similar argument to that presented in the case $i=0$.

If~(2), i.e., if $\mathfrak{P}^{(2)}$ is an echelonless second-order path, it could be the case that~(2.1) $\mathfrak{P}^{(2)}$ is an echelonless second-order path that is not head-constant, or~(2.2) $\mathfrak{P}^{(2)}$ is a head-constant echelonless second-order path that is not coherent or~(2.3) $\mathfrak{P}^{(2)}$ is a coherent head-constant echelonless second-order path.

If~(2.1), i.e., if $\mathfrak{P}^{(2)}$ is an echelonless second-order path that is not head-constant, then we let $i\in\bb{\mathfrak{P}^{(2)}}$ be the greatest index for which $\mathfrak{P}^{(2),0,i}$ is a head-constant echelonless second-order path. Then, according to Definition~\ref{DDCH}, we have that the value of the second-order Curry-Howard mapping at $\mathfrak{P}^{(2)}$ is given by
$$
\mathrm{CH}^{(2)}_{s}\left(
\mathfrak{P}^{(2)}
\right)
=
\mathrm{CH}^{(2)}_{s}\left(\mathfrak{P}^{(2),i+1,\bb{\mathfrak{P}^{(2)}}-1}\right)
\circ_{s}^{1\mathbf{T}_{\Sigma^{\boldsymbol{\mathcal{A}}^{(2)}}}(X)}
\mathrm{CH}^{(2)}_{s}\left(\mathfrak{P}^{(2),0,i}\right).
$$

Therefore,  the following chain of equalities holds
\begin{flushleft}
$\mathrm{ip}^{(2,X)@}_{s}\left(
\mathrm{CH}^{(2)}_{s}\left(
\mathfrak{P}^{(2)}
\right)\right)
$
\begin{align*}
&=
\mathrm{ip}^{(2,X)@}_{s}\left(
\mathrm{CH}^{(2)}_{s}\left(
\mathfrak{P}^{(2),i+1,\bb{\mathfrak{P}^{(2)}}-1}
\right)
\circ_{s}^{1\mathbf{T}_{\Sigma^{\boldsymbol{\mathcal{A}}^{(2)}}}(X)}
\mathrm{CH}^{(2)}_{s}\left(
\mathfrak{P}^{(2),0,i}
\right)\right)
\tag{1}
\\
&=
\mathrm{ip}^{(2,X)@}_{s}\left(
\mathrm{CH}^{(2)}_{s}\left(
\mathfrak{P}^{(2),i+1,\bb{\mathfrak{P}^{(2)}}-1}
\right)\right)
\circ_{s}^{1\mathbf{F}_{\Sigma^{\boldsymbol{\mathcal{A}}^{(2)}}}
\left(\mathbf{Pth}_{\boldsymbol{\mathcal{A}}^{(2)}}\right)
}
\\&\qquad\qquad\qquad\qquad\qquad\qquad\qquad\qquad\qquad\qquad
\mathrm{ip}^{(2,X)@}_{s}\left(
\mathrm{CH}^{(2)}_{s}\left(
\mathfrak{P}^{(2),0,i}
\right)\right)
\tag{2}
\\&=
\mathrm{ip}^{(2,X)@}_{s}\left(
\mathrm{CH}^{(2)}_{s}\left(
\mathfrak{P}^{(2),i+1,\bb{\mathfrak{P}^{(2)}}-1}
\right)\right)
\circ_{s}^{1\mathbf{Pth}_{\boldsymbol{\mathcal{A}}^{(2)}}
}
\mathrm{ip}^{(2,X)@}_{s}\left(
\mathrm{CH}^{(2)}_{s}\left(
\mathfrak{P}^{(2),0,i}
\right)\right).
\tag{3}
\end{align*}
\end{flushleft}

In the just stated chain of equalities, the first equality recovers the value of the second-order Curry-Howard mapping at $\mathfrak{P}^{(2)}$; the second equality holds since, according to Definition~\ref{DDIp}, $\mathrm{ip}^{(2,X)@}$ is a many-sorted $\Sigma^{\boldsymbol{\mathcal{A}}^{(2)}}$-homomorphism; finally, for the last equality let us note that $(\mathfrak{P}^{(2)i+1,\bb{\mathfrak{P}^{(2)}}-1},s)$ and $(\mathfrak{P}^{(2),0,i},s)$ are strictly smaller than $(\mathfrak{P}^{(2)},s)$ with respect to $\prec_{\mathbf{Pth}_{\boldsymbol{\mathcal{A}}^{(2)}}}$. Then, by the inductive hypothesis, we have that $\mathrm{ip}^{(2,X)@}_{s}(\mathrm{CH}^{(2)}_{s}(\mathfrak{P}^{(2),i+1,\bb{\mathfrak{P}^{(2)}}-1}))$ and $\mathrm{ip}^{(2,X)@}_{s}(\mathrm{CH}^{(2)}_{s}(\mathfrak{P}^{(2),0,i}))$ are second-order paths in $[\mathfrak{P}^{(2),i+1,\bb{\mathfrak{P}^{(2)}}-1}]^{}_{s}$ and $[\mathfrak{P}^{(2),0,i}]^{}_{s}$ respectively. Hence, the interpretation of the $1$-composition operation symbol $\circ^{1}_{s}$ in $\mathbf{F}_{\Sigma^{\boldsymbol{\mathcal{A}}^{(2)}}}(\mathbf{Pth}_{\boldsymbol{\mathcal{A}}^{(2)}})$ becomes that of $\mathbf{Pth}_{\boldsymbol{\mathcal{A}}^{(2)}}$. Note that the $1$-composition is defined because second-order paths in the same class have the same $([1],2)$-source and $([1],2)$-target according to Lemma~\ref{LDCH}.

Hence $\mathrm{ip}^{(1,X)@}_{s}(\mathrm{CH}^{(2)}_{s}(\mathfrak{P}^{(2)}))$ is a second-order path.

Since $\mathfrak{P}^{(2),0,i}$ is a head-constant echelonless second-order path, for a unique word $\mathbf{s}\in S^{\star}-\{\lambda\}$, and a unique operation symbol $\tau\in\Sigma^{\boldsymbol{\mathcal{A}}}_{\mathbf{s},s}$, we have that the family of first-order translations occurring in $\mathfrak{P}^{(2),0,i}$ is a family of first-order translations of type $\tau$.

By Lemma~\ref{LDCHNEchHd}, we have that $\mathrm{CH}^{(2)}_{s}(\mathfrak{P}^{(2),0,i})\in\mathcal{T}(\tau,\mathrm{T}_{\Sigma^{\boldsymbol{\mathcal{A}}^{(2)}}}(X))^{*}$, which is a subset of $[\mathrm{T}_{\Sigma^{\boldsymbol{\mathcal{A}}^{(2)}}}(X)]^{\mathsf{HdC}}_{s}$. Since $\mathrm{ip}^{(2,X)@}_{s}(\mathrm{CH}^{(2)}_{s}(\mathfrak{P}^{(2),0,i}))$ is a second-order path in $[\mathfrak{P}^{(2),0,i}]^{}_{s}$, we have, by Lemma~\ref{LDCHNEchHd}, that $\mathrm{ip}^{(2,X)@}_{s}(\mathrm{CH}^{(2)}_{s}(\mathfrak{P}^{(2),0,i}))$ is a head-constant echelonless second-order path associated to the operation symbol $\tau$ in $\Sigma^{\boldsymbol{\mathcal{A}}}_{\mathbf{s},s}$.

Since $\mathfrak{P}^{(2),i+1,\bb{\mathfrak{P}^{(2)}}-1}$ is an echelonless second-order path, we have that, for a unique word $\mathbf{s}'\in S^{\star}-\{\lambda\}$ and a unique operation symbol $\tau'\in\Sigma^{\boldsymbol{\mathcal{A}}}_{\mathbf{s}', s}$,  we have that the family of first-order translations occurring in $\mathfrak{P}^{(2),i+1,\bb{\mathfrak{P}^{(2)}}-1}$ is a family of first-order translations of type $\tau'$.

Note that $\tau\neq\tau'$ since $\mathfrak{P}^{(2),0,i}$ is the greatest subpath of $\mathfrak{P}^{(2)}$ that is head-constant. From Lemma~\ref{LDCHNEch}, we have that $\mathrm{CH}^{(2)}_{s}(\mathfrak{P}^{(2),i+1,\bb{\mathfrak{P}^{(2)}}-1})$ is a term in $\mathrm{T}_{\Sigma^{\boldsymbol{\mathcal{A}}^{(2)}}}(X)_{s}$ that does not belong to either $\eta^{(2,1)\sharp}_{s}[\mathrm{PT}_{\boldsymbol{\mathcal{A}}}]_{s}$ or $\eta^{(2,\mathcal{A}^{(2)})}[\mathcal{A}^{(2)}]^{\mathrm{pct}}_{s}$. Since $\mathrm{ip}^{(2,X)@}_{s}(\mathrm{CH}^{(2)}_{s}(\mathfrak{P}^{(2),i+1,\bb{\mathfrak{P}^{(2)}}-1}))$ is a second-order path in $[\mathfrak{P}^{(2),i+1,\bb{\mathfrak{P}^{(2)}}-1}]^{}_{s}$, we have, by Lemma~\ref{LDCHNEch}, that $\mathrm{ip}^{(2,X)@}_{s}(\mathrm{CH}^{(2)}_{s}(\mathfrak{P}^{(2),i+1,\bb{\mathfrak{P}^{(2)}}-1}))$  is an echelonless second-order path whose initial path term representative is associated to the operation symbol $\tau'$ in $\Sigma^{\boldsymbol{\mathcal{A}}}_{\mathbf{s}',s}$.

All in all, we can affirm that $\mathrm{ip}^{(2,X)@}_{s}(\mathrm{CH}^{(2)}_{s}(\mathfrak{P}^{(2)}))$ is an echelonless second-order path that is not head-constant. Moreover, $i\in\bb{\mathfrak{P}^{(2)}}$ is the greatest index for which $\mathrm{ip}^{(2,X)@}_{s}(\mathrm{CH}^{(2)}_{s}(\mathfrak{P}^{(2)}))^{0,i}$ is a head-constant echelonless second-order path. Hence, according to Definition~\ref{DDCH}, the value of the second-order Curry-Howard mapping at $\mathrm{ip}^{(2,X)@}_{s}(\mathrm{CH}^{(2)}_{s}(\mathfrak{P}^{(2)}))$ is given by 
\begin{flushleft}
$
\mathrm{CH}^{(2)}_{s}\left(
\mathrm{ip}^{(2,X)@}_{s}\left(
\mathrm{CH}^{(2)}_{s}\left(
\mathfrak{P}^{(2)}
\right)\right)\right)
$
\begin{align*}
\quad&=
\mathrm{CH}^{(2)}_{s}\left(
\left(
\mathrm{ip}^{(2,X)@}_{s}\left(
\mathrm{CH}^{(2)}_{s}\left(
\mathfrak{P}^{(2)}
\right)\right)
\right)^{i+1,\bb{
\mathrm{ip}^{(2,X)@}_{s}\left(
\mathrm{CH}^{(2)}_{s}\left(
\mathfrak{P}^{(2)}
\right)\right)
}-1}
\right)
\circ^{1\mathbf{T}_{\Sigma^{\boldsymbol{\mathcal{A}}^{(2)}}}(X)}_{s}
\\&\qquad\qquad\qquad\qquad\qquad\qquad\qquad
\mathrm{CH}^{(2)}_{s}\left(
\left(
\mathrm{ip}^{(2,X)@}_{s}\left(
\mathrm{CH}^{(2)}_{s}\left(
\mathfrak{P}^{(2)}
\right)\right)
\right)^{0,i}\right)
\tag{1}
\\&=
\mathrm{CH}^{(2)}_{s}\left(
\mathrm{ip}^{(2,X)@}_{s}\left(
\mathrm{CH}^{(2)}_{s}\left(
\mathfrak{P}^{(2),i+1,\bb{\mathfrak{P}^{(2)}}-1}
\right)\right)\right)
\circ_{s}^{1
\mathbf{T}_{\Sigma^{\boldsymbol{\mathcal{A}}^{(2)}}}(X)
}
\\&\qquad\qquad\qquad\qquad\qquad\qquad\qquad\qquad
\mathrm{CH}^{(2)}_{s}\left(
\mathrm{ip}^{(2,X)@}_{s}\left(
\mathrm{CH}^{(2)}_{s}\left(
\mathfrak{P}^{(2),0,i}
\right)\right)\right)
\tag{2}
\\&=
\mathrm{CH}^{(2)}_{s}\left(
\mathfrak{P}^{(2),i+1,\bb{\mathfrak{P}^{(2)}}-1}
\right)
\circ_{s}^{1
\mathbf{T}_{\Sigma^{\boldsymbol{\mathcal{A}}^{(2)}}}(X)
}
\mathrm{CH}^{(2)}_{s}\left(
\mathfrak{P}^{(2),0,i}
\right)
\tag{3}
\\&=
\mathrm{CH}^{(2)}_{s}\left(
\mathfrak{P}^{(2)}
\right).
\tag{4}
\end{align*}
\end{flushleft}

The first equality follows from the fact that $\mathrm{ip}^{(2,X)@}_{s}(\mathrm{CH}^{(2)}_{s}(\mathfrak{P}^{(2)}))$ is an echelonless second-order path of length strictly greater than one that is not head-constant and $i$ is the greatest index for which $\mathrm{ip}^{(2,X)@}_{s}(\mathrm{CH}^{(2)}_{s}(\mathfrak{P}^{(2)}))^{0,i}$ is a head-constant second-order path and Definition~\ref{DDCH}; the second equality follows from the fact that, from the previous discussion, $\mathrm{ip}^{(2,X)@}_{s}(\mathrm{CH}^{(2)}_{s}(\mathfrak{P}^{(2),0,i}))$ is a head-constant echelonless second-order path. Hence, from the fact that $\mathrm{ip}^{(2,X)@}_{s}(\mathrm{CH}^{(2)}_{s}(\mathfrak{P}^{(2)}))$ can be written as the $1$-composition
$$
\mathrm{ip}^{(2,X)@}_{s}\left(
\mathrm{CH}^{(2)}_{s}\left(
\mathfrak{P}^{(2),i+1,\bb{\mathfrak{P}^{(2)}}-1}
\right)\right)
\circ_{s}^{1\mathbf{Pth}_{\boldsymbol{\mathcal{A}}^{(2)}}}
\mathrm{ip}^{(2,X)@}_{s}\left(
\mathrm{CH}^{(2)}_{s}\left(
\mathfrak{P}^{(2),0,i}
\right)\right),
$$
we infer that 
\begin{align*}
&
\mathrm{ip}^{(2,X)@}_{s}\left(
\mathrm{CH}^{(2)}_{s}\left(
\mathfrak{P}^{(2)}
\right)\right)^{0,i}
=
\mathrm{ip}^{(2,X)@}_{s}\left(
\mathrm{CH}^{(2)}_{s}\left(
\mathfrak{P}^{(2),0,i}
\right)\right),  \mbox{ and} 
\\&
\mathrm{ip}^{(2,X)@}_{s}\left(
\mathrm{CH}^{(2)}_{s}\left(
\mathfrak{P}^{(2)}
\right)\right)^{i+1,\bb{\mathrm{ip}^{(2,X)@}_{s}(
\mathrm{CH}^{(2)}_{s}(
\mathfrak{P}^{(2)}
))}-1}
=
\mathrm{ip}^{(2,X)@}_{s}\left(
\mathrm{CH}^{(2)}_{s}\left(
\mathfrak{P}^{(2),i+1,\bb{\mathfrak{P}^{(2)}}-1}
\right)\right);
\end{align*}
the third equality follows from the fact that $\mathrm{ip}^{(2,X)@}_{s}(\mathrm{CH}^{(2)}_{s}(\mathfrak{P}^{(2),i+1,\bb{\mathfrak{P}^{(2)}}-1}))$ and $\mathrm{ip}^{(2,X)@}_{s}(\mathrm{CH}^{(2)}_{s}(\mathfrak{P}^{(2),0,i}))$ are second-order paths in $[\mathfrak{P}^{(2),i+1,\bb{\mathfrak{P}^{(2)}}-1}]^{}_{s}$
and $[\mathfrak{P}^{(2),0,i}]^{}_{s}$ respectively; finally, the last equality recovers the value of the second-order Curry-Howard mapping at $\mathfrak{P}^{(2)}$.

Therefore,
\begin{itemize}
\item[(i)] $\mathrm{ip}^{(2,X)@}_{s}(\mathrm{CH}^{(2)}_{s}(\mathfrak{P}^{(2)}))$ is a second-order path in $\mathrm{Pth}_{\boldsymbol{\mathcal{A}}^{(2)},s}$; and
\item[(ii)] $\mathrm{ip}^{(2,X)@}_{s}(\mathrm{CH}^{(2)}_{s}(\mathfrak{P}^{(2)}))$ belongs to $[\mathfrak{P}^{(2)}]^{}_{s}$.
\end{itemize}

The case of $\mathfrak{P}^{(2)}$ being an echelonless second-order path that is not head-constant follows.

If~(2.2), i.e., if $\mathfrak{P}^{(2)}$ is a head-constant  echelonless second-order path that is not coherent, then we let $i\in\bb{\mathfrak{P}^{(2)}}$ be the greatest index for which $\mathfrak{P}^{(2),0,i}$ is a  coherent head-constant  echelonless second-order path. Then, according to Definition~\ref{DDCH}, we have that the value of the second-order Curry-Howard mapping at $\mathfrak{P}^{(2)}$ is given by 
$$
\mathrm{CH}^{(2)}_{s}\left(
\mathfrak{P}^{(2)}
\right)
=
\mathrm{CH}^{(2)}_{s}\left(
\mathfrak{P}^{(2),i+1,\bb{\mathfrak{P}^{(2)}}-1}
\right)
\circ_{s}^{1\mathbf{T}_{\Sigma^{\boldsymbol{\mathcal{A}}^{(2)}}}(X)}
\mathrm{CH}^{(2)}_{s}\left(
\mathfrak{P}^{(2),0,i}
\right).
$$

Therefore, the following chain of equalities holds
\begin{flushleft}
$\mathrm{ip}^{(2,X)@}_{s}\left(
\mathrm{CH}^{(2)}_{s}\left(
\mathfrak{P}^{(2)}
\right)\right)
$
\begin{align*}
&=
\mathrm{ip}^{(2,X)@}_{s}\left(
\mathrm{CH}^{(2)}_{s}\left(
\mathfrak{P}^{(2),i+1,\bb{\mathfrak{P}^{(2)}}-1}
\right)
\circ_{s}^{1\mathbf{T}_{\Sigma^{\boldsymbol{\mathcal{A}}^{(2)}}}(X)}
\mathrm{CH}^{(2)}_{s}\left(
\mathfrak{P}^{(2),0,i}
\right)\right)
\tag{1}
\\
&=
\mathrm{ip}^{(2,X)@}_{s}\left(
\mathrm{CH}^{(2)}_{s}\left(
\mathfrak{P}^{(2),i+1,\bb{\mathfrak{P}^{(2)}}-1}
\right)\right)
\circ_{s}^{1\mathbf{F}_{\Sigma^{\boldsymbol{\mathcal{A}}^{(2)}}}
\left(\mathbf{Pth}_{\boldsymbol{\mathcal{A}}^{(2)}}\right)
}
\\&\qquad\qquad\qquad\qquad\qquad\qquad\qquad\qquad\qquad\qquad
\mathrm{ip}^{(2,X)@}_{s}\left(
\mathrm{CH}^{(2)}_{s}\left(
\mathfrak{P}^{(2),0,i}
\right)\right)
\tag{2}
\\&=
\mathrm{ip}^{(2,X)@}_{s}\left(
\mathrm{CH}^{(2)}_{s}\left(
\mathfrak{P}^{(2),i+1,\bb{\mathfrak{P}^{(2)}}-1}
\right)\right)
\circ_{s}^{1\mathbf{Pth}_{\boldsymbol{\mathcal{A}}^{(2)}}
}
\mathrm{ip}^{(2,X)@}_{s}\left(
\mathrm{CH}^{(2)}_{s}\left(
\mathfrak{P}^{(2),0,i}
\right)\right).
\tag{3}
\end{align*}
\end{flushleft}

In the just stated chain of equalities, the first equality recovers the value of the second-order Curry-Howard mapping at $\mathfrak{P}^{(2)}$; the second equality holds since, according to Definition~\ref{DDIp}, $\mathrm{ip}^{(2,X)@}$ is a $\Sigma^{\boldsymbol{\mathcal{A}}^{(2)}}$-homomorphism; finally, for the last equality let us note that $(\mathfrak{P}^{(2)i+1,\bb{\mathfrak{P}^{(2)}}-1},s)$ and $(\mathfrak{P}^{(2),0,i},s)$ are strictly smaller than $(\mathfrak{P}^{(2)},s)$ with respect to $\prec_{\mathbf{Pth}_{\boldsymbol{\mathcal{A}}^{(2)}}}$. Then, by the inductive hypothesis, we have that $\mathrm{ip}^{(2,X)@}_{s}(\mathrm{CH}^{(2)}_{s}(\mathfrak{P}^{(2),i+1,\bb{\mathfrak{P}^{(2)}}-1}))$ and $\mathrm{ip}^{(2,X)@}_{s}(\mathrm{CH}^{(2)}_{s}(\mathfrak{P}^{(2),0,i}))$ are second-order paths in $[\mathfrak{P}^{(2),i+1,\bb{\mathfrak{P}^{(2)}}-1}]^{}_{s}$ and $[\mathfrak{P}^{(2),0,i}]^{}_{s}$ respectively. Hence, the interpretation of the $1$-composition operation symbol $\circ^{1}_{s}$ in $\mathbf{F}_{\Sigma^{\boldsymbol{\mathcal{A}}^{(2)}}}(\mathbf{Pth}_{\boldsymbol{\mathcal{A}}^{(2)}})$ becomes that of $\mathbf{Pth}_{\boldsymbol{\mathcal{A}}^{(2)}}$. Note that the $1$-composition is defined because second-order paths in the same class have the same $([1],2)$-source and $([1],2)$-target according to Lemma~\ref{LDCH}. 

Hence $\mathrm{ip}^{(2,X)@}_{s}(\mathrm{CH}^{(2)}_{s}(\mathfrak{P}^{(2)}))$ is a proper second-order path.

Let us recall that $\mathfrak{P}^{(2)}$ can be written as the $1$-composition 
$
\mathfrak{P}^{(2),i+1,\bb{\mathfrak{P}^{(2)}}-1}
\circ^{1\mathbf{Pth}_{\boldsymbol{\mathcal{A}}^{(2)}}}_{s}
\mathfrak{P}^{(2),0,i}
$ 
and $\mathfrak{P}^{(2)}$ is a head-constant echelonless second-order path satisfying that $i\in\bb{\mathfrak{P}^{(2)}}$ is the greatest index for which $\mathfrak{P}^{(2),0,i}$ is coherent. Moreover, we have by induction that $\mathrm{ip}^{(2,X)@}_{s}(\mathrm{CH}^{(2)}_{s}(\mathfrak{P}^{(2),i+1,\bb{\mathfrak{P}^{(2)}}-1}))$ and $\mathrm{ip}^{(2,X)@}_{s}(\mathrm{CH}^{(2)}_{s}(\mathfrak{P}^{(2),0,i}))$ are second-order paths in $[\mathfrak{P}^{(2),i+1,\bb{\mathfrak{P}^{(2)}}-1}]^{}_{s}$ and $[\mathfrak{P}^{(2),0,i}]^{}_{s}$ respectively. Therefore, since $\mathrm{ip}^{(2,X)@}_{s}(\mathrm{CH}^{(2)}_{s}(\mathfrak{P}^{(2)}))$, by the above equality, decomposes as the $1$-composition
$$
\mathrm{ip}^{(2,X)@}_{s}\left(
\mathrm{CH}^{(2)}_{s}\left(
\mathfrak{P}^{(2),i+1,\bb{\mathfrak{P}^{(2)}}-1}
\right)\right)
\circ_{s}^{1\mathbf{Pth}_{\boldsymbol{\mathcal{A}}^{(2)}}
}
\mathrm{ip}^{(2,X)@}_{s}\left(
\mathrm{CH}^{(2)}_{s}\left(
\mathfrak{P}^{(2),0,i}
\right)\right),
$$
we conclude, in virtue of Lemma~\ref{LTech}, that $\mathrm{ip}^{(2,X)@}_{s}(\mathrm{CH}^{(2)}_{s}(\mathfrak{P}^{(2)}))$ is a head-constant echelonless second-order path that is not coherent. Moreover $i\in\bb{\mathfrak{P}^{(2)}}$ is the greatest index for which $\mathrm{ip}^{(2,X)@}_{s}(\mathrm{CH}^{(2)}_{s}(\mathfrak{P}^{(2)}))^{0,i}$ is coherent.

Hence, according to Definition~\ref{DDCH}, the value of the second-order Curry-Howard mapping at $\mathrm{ip}^{(2,X)@}_{s}(\mathrm{CH}^{(2)}_{s}(\mathfrak{P}^{(2)}))$ is given by
\begin{flushleft}
$
\mathrm{CH}^{(2)}_{s}\left(
\mathrm{ip}^{(2,X)@}_{s}\left(
\mathrm{CH}^{(2)}_{s}\left(
\mathfrak{P}^{(2)}
\right)\right)\right)
$
\begin{align*}
\quad&=
\mathrm{CH}^{(2)}_{s}\left(
\left(
\mathrm{ip}^{(2,X)@}_{s}\left(
\mathrm{CH}^{(2)}_{s}\left(
\mathfrak{P}^{(2)}
\right)\right)
\right)^{i+1,\bb{
\mathrm{ip}^{(2,X)@}_{s}(
\mathrm{CH}^{(2)}_{s}(
\mathfrak{P}^{(2)}
))
}-1}
\right)
\circ^{1\mathbf{T}_{\Sigma^{\boldsymbol{\mathcal{A}}^{(2)}}}(X)}_{s}
\\&\qquad\qquad\qquad\qquad\qquad\qquad\qquad
\mathrm{CH}^{(2)}_{s}\left(
\left(
\mathrm{ip}^{(2,X)@}_{s}\left(
\mathrm{CH}^{(2)}_{s}\left(
\mathfrak{P}^{(2)}
\right)\right)
\right)^{0,i}\right)
\tag{1}
\\&=
\mathrm{CH}^{(2)}_{s}\left(
\mathrm{ip}^{(2,X)@}_{s}\left(
\mathrm{CH}^{(2)}_{s}\left(
\mathfrak{P}^{(2),i+1,\bb{\mathfrak{P}^{(2)}}-1}
\right)\right)\right)
\circ_{s}^{1
\mathbf{T}_{\Sigma^{\boldsymbol{\mathcal{A}}^{(2)}}}(X)
}
\\&\qquad\qquad\qquad\qquad\qquad\qquad\qquad\qquad
\mathrm{CH}^{(2)}_{s}\left(
\mathrm{ip}^{(2,X)@}_{s}\left(
\mathrm{CH}^{(2)}_{s}\left(
\mathfrak{P}^{(2),0,i}
\right)\right)\right)
\tag{2}
\\&=
\mathrm{CH}^{(2)}_{s}\left(
\mathfrak{P}^{(2),i+1,\bb{\mathfrak{P}^{(2)}}-1}
\right)
\circ_{s}^{1
\mathbf{T}_{\Sigma^{\boldsymbol{\mathcal{A}}^{(2)}}}(X)
}
\mathrm{CH}^{(2)}_{s}\left(
\mathfrak{P}^{(2),0,i}
\right)
\tag{3}
\\&=
\mathrm{CH}^{(2)}_{s}\left(
\mathfrak{P}^{(2)}
\right).
\tag{4}
\end{align*}
\end{flushleft}

The first equality follows from the fact that $\mathrm{ip}^{(2,X)@}_{s}(\mathrm{CH}^{(2)}_{s}(\mathfrak{P}^{(2)}))$ is a head-constant echelonless second-order path that is not head-constant and $i$ is the greatest index for which $\mathrm{ip}^{(2,X)@}_{s}(\mathrm{CH}^{(2)}_{s}(\mathfrak{P}^{(2)}))^{0,i}$ is a head-constant, coherent second-order path and Definition~\ref{DDCH}; the second equality follows from the fact that, from the previous discussion, $\mathrm{ip}^{(2,X)@}_{s}(\mathrm{CH}^{(2)}_{s}(\mathfrak{P}^{(2),0,i}))$ is a coherent head-constant echelonless second-order path. Hence, from the fact that $\mathrm{ip}^{(2,X)@}_{s}(\mathrm{CH}^{(2)}_{s}(\mathfrak{P}^{(2)}))$ can be written as the $1$-composition
$$
\mathrm{ip}^{(2,X)@}_{s}\left(
\mathrm{CH}^{(2)}_{s}\left(
\mathfrak{P}^{(2),i+1,\bb{\mathfrak{P}^{(2)}}-1}
\right)\right)
\circ_{s}^{1\mathbf{Pth}_{\boldsymbol{\mathcal{A}}^{(2)}}}
\mathrm{ip}^{(2,X)@}_{s}\left(
\mathrm{CH}^{(2)}_{s}\left(
\mathfrak{P}^{(2),0,i}
\right)\right),
$$
we infer that 
\begin{align*}
&
\mathrm{ip}^{(2,X)@}_{s}\left(
\mathrm{CH}^{(2)}_{s}\left(
\mathfrak{P}^{(2)}
\right)\right)^{0,i}
=
\mathrm{ip}^{(2,X)@}_{s}\left(
\mathrm{CH}^{(2)}_{s}\left(
\mathfrak{P}^{(2),0,i}
\right)\right),  \mbox{ and} 
\\&
\mathrm{ip}^{(2,X)@}_{s}\left(
\mathrm{CH}^{(2)}_{s}\left(
\mathfrak{P}^{(2)}
\right)\right)^{i+1,\bb{\mathrm{ip}^{(2,X)@}_{s}(
\mathrm{CH}^{(2)}_{s}(
\mathfrak{P}^{(2)}
))}-1}
=
\mathrm{ip}^{(2,X)@}_{s}\left(
\mathrm{CH}^{(2)}_{s}\left(
\mathfrak{P}^{(2),i+1,\bb{\mathfrak{P}^{(2)}}-1}
\right)\right);
\end{align*}
the third equality follows from the fact that $\mathrm{ip}^{(2,X)@}_{s}(\mathrm{CH}^{(2)}_{s}(\mathfrak{P}^{(2),i+1,\bb{\mathfrak{P}^{(2)}}-1}))$ and $\mathrm{ip}^{(2,X)@}_{s}(\mathrm{CH}^{(2)}_{s}(\mathfrak{P}^{(2),0,i}))$ are second-order paths in $[\mathfrak{P}^{(2),i+1,\bb{\mathfrak{P}^{(2)}}-1}]^{}_{s}$
and $[\mathfrak{P}^{(2),0,i}]^{}_{s}$, respectively; finally, the last equality recovers the value of the second-order Curry-Howard mapping at $\mathfrak{P}^{(2)}$.

Therefore,
\begin{itemize}
\item[(i)] $\mathrm{ip}^{(2,X)@}_{s}(\mathrm{CH}^{(2)}_{s}(\mathfrak{P}^{(2)}))$ is a second-order path in $\mathrm{Pth}_{\boldsymbol{\mathcal{A}}^{(2)},s}$; and
\item[(ii)] $\mathrm{ip}^{(2,X)@}_{s}(\mathrm{CH}^{(2)}_{s}(\mathfrak{P}^{(2)}))$ belongs to $[\mathfrak{P}^{(2)}]^{}_{s}$.
\end{itemize}

The case of $\mathfrak{P}^{(2)}$ being a head-constant echelonless second-order path that is not coherent follows.

If~(2.3), i.e., if $\mathfrak{P}^{(2)}$ is a coherent  head-constant echelonless second-order path, then, by Definition~\ref{DDHeadCt}, there exists a unique word $\mathbf{s}\in S^{\star}-\{\lambda\}$ and a unique operation symbol $\tau$ in $\Sigma^{\boldsymbol{\mathcal{A}}}_{\mathbf{s},s}$ associated to $\mathfrak{P}^{(2)}$. Let us recall from Definition~\ref{DDPth} that the operations of $0$-source and $0$-target are forbidden. Let $(\mathfrak{P}^{(2)}_{j})_{j\in\bb{\mathbf{s}}}$ be the family of second-order paths we can extract from $\mathfrak{P}^{(2)}$ in virtue of Lemma~\ref{LDPthExtract}. Then, according to Definition~\ref{DDCH}, we have that the value of the second-order Curry-Howard mapping at $\mathfrak{P}^{(2)}$ is given by
$$
\mathrm{CH}^{(2)}_{s}\left(
\mathfrak{P}^{(2)}
\right)
=
\tau^{\mathbf{T}_{\Sigma^{\boldsymbol{\mathcal{A}}^{(2)}}}(X)}
\left(\left(\mathrm{CH}^{(2)}_{s_{j}}\left(
\mathfrak{P}^{(2)}_{j}
\right)\right)_{j\in\bb{\mathbf{s}}}
\right).
$$

Therefore, the following chain of equalities holds
\begin{flushleft}
$\mathrm{ip}^{(2,X)@}_{s}\left(
\mathrm{CH}^{(2)}_{s}\left(
\mathfrak{P}^{(2)}
\right)\right)$
\allowdisplaybreaks
\begin{align*}
\qquad&=
\mathrm{ip}^{(2,X)@}_{s}\left(
\tau^{\mathbf{T}_{\Sigma^{\boldsymbol{\mathcal{A}}^{(2)}}}(X)}
\left(\left(\mathrm{CH}^{(2)}_{s_{j}}\left(
\mathfrak{P}^{(2)}_{j}
\right)\right)_{j\in\bb{\mathbf{s}}}\right)
\right)
\tag{1}
\\&=
\tau^{\mathbf{F}_{\Sigma^{\boldsymbol{\mathcal{A}}^{(2)}}}
\left(\mathbf{Pth}_{\boldsymbol{\mathcal{A}}^{(2)}}
\right)
}
\left(\left(
\mathrm{ip}^{(2,X)@}_{s_{j}}\left(
\mathrm{CH}^{(2)}_{s_{j}}\left(
\mathfrak{P}^{(2)}_{j}
\right)\right)
\right)_{j\in\bb{\mathbf{s}}}\right)
\tag{2}
\\&=
\tau^{
\mathbf{Pth}_{\boldsymbol{\mathcal{A}}^{(2)}}
}
\left(\left(
\mathrm{ip}^{(2,X)@}_{s_{j}}\left(
\mathrm{CH}^{(2)}_{s_{j}}\left(
\mathfrak{P}^{(2)}_{j}
\right)\right)
\right)_{j\in\bb{\mathbf{s}}}\right).
\tag{3}
\end{align*}
\end{flushleft}

In the just stated chain of equalities, the first equality recovers the value of the second-order Curry-Howard mapping at $\mathfrak{P}^{(2)}$; the second equality holds since, according to Definition~\ref{DDIp}, $\mathrm{ip}^{(2,X)@}$ is a $\Sigma^{\boldsymbol{\mathcal{A}}^{(2)}}$-homomorphism; finally, for the last equality let us note that, for every $j\in\bb{\mathbf{s}}$, $(\mathfrak{P}^{(2)}_{j},s_{j})$ is strictly smaller than $(\mathfrak{P}^{(2)},s)$ with respect to $\prec_{\mathbf{Pth}_{\boldsymbol{\mathcal{A}}^{(2)}}}$. Then, by the inductive hypothesis, we have that $\mathrm{ip}^{(2,X)@}_{s_{j}}(
\mathrm{CH}^{(2)}_{s_{j}}(
\mathfrak{P}^{(2)}_{j}
))$ is a second-order path in $[\mathfrak{P}^{(2)}_{j}]^{}_{s_{j}}$. Hence, the interpretation of the operation symbol $\tau$ in $\mathbf{F}_{\Sigma^{\boldsymbol{\mathcal{A}}^{(2)}}}(
\mathbf{Pth}_{\boldsymbol{\mathcal{A}}^{(2)}}
)
$ becomes that of $\mathbf{Pth}_{\boldsymbol{\mathcal{A}}^{(2)}}$. Note that this operation is well-defined because second-order paths in the same class have the same $([1],2)$-source and $([1],2)$-target according to Lemma~\ref{LDCH}. 

Hence $\mathrm{ip}^{(2,X)@}_{s}(\mathrm{CH}^{(2)}_{s}(\mathfrak{P}^{(2)}))$ is a second-order path.

Since $\mathfrak{P}^{(2)}$ is a second-order path of length at least one,  there exists at least one index $j\in\bb{\mathbf{s}}$ for which $\mathfrak{P}^{(2)}_{j}$ has length at least one. Then, since $(\mathfrak{P}^{(2)}_{j},\mathrm{ip}^{(2,X)@}_{s_{j}}(\mathrm{CH}^{(2)}_{s_{j}}(\mathfrak{P}^{(2)}_{j})))$ is a pair in $\mathrm{Ker}(\mathrm{CH}^{(2)})_{s_{j}}$, we have, by Lemma~\ref{LDCH}, that $\mathrm{ip}^{(2,X)@}_{s_{j}}(\mathrm{CH}^{(2)}_{s_{j}}(\mathfrak{P}^{(2)}_{j}))$ is a second-order path of length at least one. Then, since 
$$
\mathrm{ip}^{(2,X)@}_{s}\left(
\mathrm{CH}^{(2)}_{s}\left(
\mathfrak{P}^{(2)}
\right)\right)=
\tau^{
\mathbf{Pth}_{\boldsymbol{\mathcal{A}}^{(2)}}
}
\left(\left(
\mathrm{ip}^{(2,X)@}_{s_{j}}\left(
\mathrm{CH}^{(2)}_{s_{j}}\left(
\mathfrak{P}^{(2)}_{j}
\right)\right)
\right)_{j\in\bb{\mathbf{s}}}\right),$$
we have, by Corollary~\ref{CDPthWB}, that $\mathrm{ip}^{(2,X)@}_{s}(
\mathrm{CH}^{(2)}_{s}(
\mathfrak{P}^{(2)}
))$ is a coherent head-constant echelonless second-order path. Moreover, in virtue of Proposition~\ref{PDRecov}, the second-order path extraction procedure from Lemma~\ref{LDPthExtract} applied to $\mathrm{ip}^{(2,X)@}_{s}(
\mathrm{CH}^{(2)}_{s}(
\mathfrak{P}^{(2)}
))$ retrieves the family of second-order paths $(
\mathrm{ip}^{(2,X)@}_{s_{j}}(
\mathrm{CH}^{(2)}_{s_{j}}(
\mathfrak{P}^{(2)}_{j}
))
)_{j\in\bb{\mathbf{s}}}$.

Hence, according to Definition~\ref{DDCH}, the value of the second-order Curry-Howard mapping at $\mathrm{ip}^{(2,X)@}_{s}(
\mathrm{CH}^{(2)}_{s}(
\mathfrak{P}^{(2)}
))$ is given by
\begin{flushleft}
$\mathrm{CH}^{(2)}_{s}\left(
\mathrm{ip}^{(2,X)@}_{s}\left(
\mathrm{CH}^{(2)}_{s}\left(
\mathfrak{P}^{(2)}
\right)\right)\right)$
\allowdisplaybreaks
\begin{align*}
\quad
&=
\tau^{\mathbf{T}_{\Sigma^{\boldsymbol{\mathcal{A}}^{(2)}}}(X)}
\left(\left(\mathrm{CH}^{(2)}_{s_{j}}\left(
\mathrm{ip}^{(2,X)@}_{s_{j}}\left(
\mathrm{CH}^{(2)}_{s_{j}}\left(
\mathfrak{P}^{(2)}_{j}
\right)\right)\right)
\right)_{j\in\bb{\mathbf{s}}}
\right)
\tag{1}
\\&=
\tau^{\mathbf{T}_{\Sigma^{\boldsymbol{\mathcal{A}}^{(2)}}}(X)}
\left(\left(\mathrm{CH}^{(2)}_{s_{j}}\left(
\mathfrak{P}^{(2)}_{j}
\right)
\right)_{j\in\bb{\mathbf{s}}}
\right)
\tag{2}
\\&=
\mathrm{CH}^{(2)}_{s}\left(
\mathfrak{P}^{(2)}
\right).
\tag{3}
\end{align*}
\end{flushleft}

The first equality follows from the fact that $\mathrm{ip}^{(2,X)@}_{s}(
\mathrm{CH}^{(2)}_{s}(
\mathfrak{P}^{(2)}
))$ is a coherent head-constant echelonless second-order path and Definition~\ref{DDCH}; the second equality follows from the fact that, from the previous discussion, for every $j\in\bb{\mathbf{s}}$, $\mathrm{ip}^{(2,X)@}_{s_{j}}(\mathrm{CH}^{(2)}_{s_{j}}(\mathfrak{P}^{(2)}_{j}))$ is a second-order path in $[\mathfrak{P}^{(2)}_{j}]^{}_{s_{j}}$; finally, the last equality recovers the value of the second-order Curry-Howard mapping at $\mathfrak{P}^{(2)}$.

Therefore,
\begin{itemize}
\item[(i)] $\mathrm{ip}^{(2,X)@}_{s}(\mathrm{CH}^{(2)}_{s}(\mathfrak{P}^{(2)}))$ is a second-order path in $\mathrm{Pth}_{\boldsymbol{\mathcal{A}}^{(2)},s}$; and
\item[(ii)] $\mathrm{ip}^{(2,X)@}_{s}(\mathrm{CH}^{(2)}_{s}(\mathfrak{P}^{(2)}))$ belongs to $[\mathfrak{P}^{(2)}]^{}_{s}$.
\end{itemize}

The case of $\mathfrak{P}^{(2)}$ being a coherent head-constant echelonless second-order path follows.

This finishes the proof.
\end{proof}

The following corollary states that $(2,[1])$-identity second-order paths are always normalized.  

\begin{restatable}{corollary}{CDIpDUId}
\label{CDIpDUId}
Let $s$ be a sort in $S$, and $\mathfrak{P}^{(2)}$ a second-order path in $\mathrm{Pth}_{\boldsymbol{\mathcal{A}}^{(2)},s}$. If $\mathfrak{P}^{(2)}$ is a  $(2,[1])$-identity second-order path, then $\mathrm{ip}^{(2,X)@}_{s}(\mathrm{CH}^{(2)}_{s}(\mathfrak{P}^{(2)}))=\mathfrak{P}^{(2)}$.
\end{restatable}
\begin{proof}
It follows from Proposition~\ref{PDIpDCH} and Corollary~\ref{CDCHUId}.
\end{proof}

We present a series of results that will be useful later on. The first result states that the composition $\mathrm{ip}^{(2,X)@}\circ \mathrm{CH}^{(2)}$ is a $\Sigma$-homomorphism.

\begin{restatable}{lemma}{LDIpDCHSigma}
\label{LDIpDCHSigma} Let $(\mathbf{s}, s)$ be a pair in  $S^{\star}\times S$, let $\sigma$ be an operation symbol in $\Sigma_{\mathbf{s},s}$ and let $(\mathfrak{P}^{(2)}_{j})_{j\in\bb{\mathbf{s}}}$ be a family of second-order paths in $\mathrm{Pth}_{\boldsymbol{\mathcal{A}}^{(2)},\mathbf{s}}$
then the following equality holds
\allowdisplaybreaks
\begin{multline*}
\mathrm{ip}^{(2,X)@}_{s}\left(
\mathrm{CH}^{(2)}_{s}\left(
\sigma^{\mathbf{Pth}_{\boldsymbol{\mathcal{A}}^{(2)}}}\left(
\left(
\mathfrak{P}^{(2)}_{j}
\right)_{j\in\bb{\mathbf{s}}}
\right)\right)\right)\\=
\sigma^{\mathbf{Pth}_{\boldsymbol{\mathcal{A}}^{(2)}}}\left(
\left(
\mathrm{ip}^{(2,X)@}_{s_{j}}\left(
\mathrm{CH}^{(2)}_{s_{j}}\left(
\mathfrak{P}^{(2)}_{j}
\right)\right)
\right)_{j\in\bb{\mathbf{s}}}
\right).
\end{multline*}
\end{restatable}
\begin{proof}
It follows from the fact that $\mathrm{CH}^{(2)}$ and $\mathrm{ip}^{(2,X)@}$ are $\Sigma$-homomorphisms according to, respectively, Proposition~\ref{PDCHHom} and Definition~\ref{DDIp}.
\end{proof}

We next prove that the normalized form of a rewrite rule is the rewrite rule itself.

\begin{restatable}{lemma}{LDIpDCHEch}
\label{LDIpDCHEch} Let $s$ be a sort in $S$ and let $\mathfrak{p}$ be a rewrite rule in $\mathcal{A}_{s}$ then the following equality holds
\[
\mathrm{ip}^{(2,X)@}_{s}\left(
\mathrm{CH}^{(2)}_{s}\left(
\mathfrak{p}^{\mathbf{Pth}_{\boldsymbol{\mathcal{A}}^{(2)}}}
\right)\right)
=
\mathfrak{p}^{\mathbf{Pth}_{\boldsymbol{\mathcal{A}}^{(2)}}}.
\]
\end{restatable}
\begin{proof}
The following chain of equalities holds.
\allowdisplaybreaks
\begin{align*}
\mathrm{ip}^{(2,X)@}_{s}\left(
\mathrm{CH}^{(2)}_{s}\left(
\mathfrak{p}^{\mathbf{Pth}_{\boldsymbol{\mathcal{A}}^{(2)}}}
\right)\right)&=
\mathrm{ip}^{(2,X)@}_{s}\left(
\mathrm{CH}^{(2)}_{s}\left(
\mathrm{ech}^{(2,\mathcal{A})}_{s}\left(
\mathfrak{p}
\right)
\right)\right)
\tag{1}
\\&=
\mathrm{ip}^{(2,X)@}_{s}\left(
\eta^{(2,\mathcal{A})}_{s}\left(
\mathfrak{p}
\right)\right)
\tag{2}
\\&=
\mathrm{ech}^{(2,\mathcal{A})}_{s}\left(
\mathfrak{p}
\right)
\tag{3}
\\&=
\mathfrak{p}^{\mathbf{Pth}_{\boldsymbol{\mathcal{A}}^{(2)}}}.
\tag{4}
\end{align*}

In the just stated chain of equalities, the first equality unravels the interpretation of the constant operation symbol $\mathfrak{p}$ in the many-sorted partial $\Sigma^{\boldsymbol{\mathcal{A}}^{(2)}}$-algebra $\mathbf{Pth}_{\boldsymbol{\mathcal{A}}^{(2)}}$, according to Proposition~\ref{PDPthCatAlg}; the second equality follows from Proposition~\ref{PDCHA}; the third equality follows from Proposition~\ref{PDIpEch}; finally, the last equality recovers the interpretation of the constant operation symbol $\mathfrak{p}$ in the many-sorted partial $\Sigma^{\boldsymbol{\mathcal{A}}^{(2)}}$-algebra $\mathbf{Pth}_{\boldsymbol{\mathcal{A}}^{(2)}}$, according to Proposition~\ref{PDPthCatAlg}.

This completes the proof.
\end{proof}

We next prove that the normalized form of a $0$-source of a second-order path is the $0$-source of the respective normalized form. The same also holds for the $0$-target.

\begin{restatable}{lemma}{LDIpDCHScZTgZ}
\label{LDIpDCHScZTgZ} Let $s$ be a sort in $S$ and let $\mathfrak{P}^{(2)}$ be a second-order path in $\mathrm{Pth}_{\boldsymbol{\mathcal{A}}^{(2)},s}$
then the following equality holds
\allowdisplaybreaks
\begin{align*}
\mathrm{ip}^{(2,X)@}_{s}\left(
\mathrm{CH}^{(2)}_{s}\left(
\mathrm{sc}^{0\mathbf{Pth}_{\boldsymbol{\mathcal{A}}^{(2)}}}_{s}\left(
\mathfrak{P}^{(2)}
\right)\right)\right)
&=
\mathrm{sc}^{0\mathbf{Pth}_{\boldsymbol{\mathcal{A}}^{(2)}}}_{s}\left(
\mathrm{ip}^{(2,X)@}_{s}\left(
\mathrm{CH}^{(2)}_{s}\left(
\mathfrak{P}^{(2)}
\right)\right)\right);
\\
\mathrm{ip}^{(2,X)@}_{s}\left(
\mathrm{CH}^{(2)}_{s}\left(
\mathrm{tg}^{0\mathbf{Pth}_{\boldsymbol{\mathcal{A}}^{(2)}}}_{s}\left(
\mathfrak{P}^{(2)}
\right)\right)\right)
&=
\mathrm{tg}^{0\mathbf{Pth}_{\boldsymbol{\mathcal{A}}^{(2)}}}_{s}\left(
\mathrm{ip}^{(2,X)@}_{s}\left(
\mathrm{CH}^{(2)}_{s}\left(
\mathfrak{P}^{(2)}
\right)\right)\right).
\end{align*}
\end{restatable}
\begin{proof}
We present the proof only for the first statement. The other one can be handled similarly.

The following chain of equalities holds
\begin{flushleft}
$\mathrm{ip}^{(2,X)@}_{s}\left(
\mathrm{CH}^{(2)}_{s}\left(
\mathrm{sc}^{0\mathbf{Pth}_{\boldsymbol{\mathcal{A}}^{(2)}}}_{s}\left(
\mathfrak{P}^{(2)}
\right)\right)\right)$
\allowdisplaybreaks
\begin{align*}
\qquad&=
\mathrm{ip}^{(2,X)@}_{s}\left(
\mathrm{CH}^{(2)}_{s}\left(
\mathrm{ip}^{(2,0)\sharp}_{s}\left(
\mathrm{sc}^{(0,2)}_{s}\left(
\mathfrak{P}^{(2)}
\right)\right)\right)\right)
\tag{1}
\\&=
\mathrm{ip}^{(2,X)@}_{s}\left(
\eta^{(2,0)}_{s}\left(
\mathrm{sc}^{(0,2)}_{s}\left(
\mathfrak{P}^{(2)}
\right)\right)\right)
\tag{2}
\\&=
\mathrm{ip}^{(2,0)\sharp}_{s}\left(
\mathrm{sc}^{(0,2)}_{s}\left(
\mathfrak{P}^{(2)}
\right)\right)
\tag{3}
\\&=
\mathrm{ip}^{(2,0)\sharp}_{s}\left(
\mathrm{sc}^{(0,2)}_{s}\left(
\mathrm{ip}^{(2,X)@}_{s}\left(
\mathrm{CH}^{(2)}_{s}\left(
\mathfrak{P}^{(2)}
\right)\right)
\right)\right)
\tag{4}
\\&=
\mathrm{sc}^{0\mathbf{Pth}_{\boldsymbol{\mathcal{A}}^{(2)}}}_{s}\left(
\mathrm{ip}^{(2,X)@}_{s}\left(
\mathrm{CH}^{(2)}_{s}\left(
\mathfrak{P}^{(2)}
\right)\right)
\right).
\tag{5}
\end{align*}
\end{flushleft}

In the just stated chain of equalities, the first equality unravels the description of the operation symbol $\mathrm{sc}^{0}_{s}$ in the many-sorted partial $\Sigma^{\boldsymbol{\mathcal{A}}^{(2)}}$-algebra $\mathbf{Pth}_{\boldsymbol{\mathcal{A}}^{(2)}}$, according to Proposition~\ref{PDPthCatAlg}; the second equality follows from Proposition~\ref{PDCHDZId}; the third equality follows from Proposition~\ref{PDIpDZ}; the fourth equality follows from Proposition~\ref{PDIpDCH} and Corollary~\ref{CDCH}; finally, the last equality recovers the description of the operation symbol $\mathrm{sc}^{0}_{s}$ in the many-sorted partial $\Sigma^{\boldsymbol{\mathcal{A}}^{(2)}}$-algebra $\mathbf{Pth}_{\boldsymbol{\mathcal{A}}^{(2)}}$, according to Proposition~\ref{PDPthCatAlg}.

This completes the proof.
\end{proof}

We next show that the normalized form of a $0$-composition of second-order paths is the $0$-composition of the respective normalized forms.

\begin{restatable}{lemma}{LDIpDCHCompZ}
\label{LDIpDCHCompZ}
Let $s$ be a sort in $S$ and let $\mathfrak{P}^{(2)}, \mathfrak{Q}^{(2)}$ be second-order paths in $\mathrm{Pth}_{\boldsymbol{\mathcal{A}}^{(2)},s}$ satisfying that 
\[
\mathrm{sc}^{(0,2)}_{s}\left(\mathfrak{Q}^{(2)}\right)
=
\mathrm{tg}^{(0,2)}_{s}\left(\mathfrak{P}^{(2)}\right),
\]
then the following equality holds
\allowdisplaybreaks
\begin{multline*}
\mathrm{ip}^{(2,X)@}_{s}\left(
\mathrm{CH}^{(2)}_{s}\left(
\mathfrak{Q}^{(2)}
\circ^{0\mathbf{Pth}_{\boldsymbol{\mathcal{A}}^{(2)}}}_{s}
\mathfrak{P}^{(2)}
\right)\right)
\\
=
\mathrm{ip}^{(2,X)@}_{s}\left(
\mathrm{CH}^{(2)}_{s}\left(
\mathfrak{Q}^{(2)}
\right)\right)
\circ^{0\mathbf{Pth}_{\boldsymbol{\mathcal{A}}^{(2)}}}_{s}
\mathrm{ip}^{(2,X)@}_{s}\left(
\mathrm{CH}^{(2)}_{s}\left(
\mathfrak{P}^{(2)}
\right)\right).
\end{multline*}
\end{restatable}
\begin{proof}
We will consider the following two cases. Either (1) both $\mathfrak{Q}^{(2)}$ and $\mathfrak{P}^{(2)}$ are $(2,[1])$-identity second-order paths, or (2) either $\mathfrak{Q}^{(2)}$ or $\mathfrak{P}^{(2)}$ is not a $(2,[1])$-identity second-order path.

For the first case, i.e., when both $\mathfrak{Q}^{(2)}$ and $\mathfrak{P}^{(2)}$ are $(2,[1])$-identity second-order paths, then there exists second-order path term classes $[Q]_{s}$ and $[P]_{s}$ in $[\mathbf{PT}_{\boldsymbol{\mathcal{A}}}]_{s}$ satisfying that 
\begin{align*}
\mathfrak{P}^{(2)}&=\mathrm{ip}^{(2,[1])\sharp}_{s}\left([P]_{s}\right);
&
\mathfrak{Q}^{(2)}&=\mathrm{ip}^{(2,[1])\sharp}_{s}\left([Q]_{s}\right).
\end{align*}

Then the following chain of equalities holds
\begin{flushleft}
$\mathrm{ip}^{(2,X)@}_{s}\left(
\mathrm{CH}^{(2)}_{s}\left(
\mathfrak{Q}^{(2)}
\circ^{0\mathbf{Pth}_{\boldsymbol{\mathcal{A}}^{(2)}}}_{s}
\mathfrak{P}^{(2)}
\right)\right)$
\allowdisplaybreaks
\begin{align*}
&=
\mathrm{ip}^{(2,X)@}_{s}\left(
\mathrm{CH}^{(2)}_{s}\left(
\mathrm{ip}^{(2,[1])\sharp}_{s}\left([Q]_{s}\right)
\circ^{0\mathbf{Pth}_{\boldsymbol{\mathcal{A}}^{(2)}}}_{s}
\mathrm{ip}^{(2,[1])\sharp}_{s}\left([P]_{s}\right)
\right)\right)
\tag{1}
\\&=
\mathrm{ip}^{(2,X)@}_{s}\left(
\mathrm{CH}^{(2)}_{s}\left(
\mathrm{ip}^{(2,[1])\sharp}_{s}\left([Q\circ^{0\mathbf{PT}_{\boldsymbol{\mathcal{A}}}}_{s} P]_{s}\right)
\right)\right)
\tag{2}
\\&=
\mathrm{ip}^{(2,[1])\sharp}_{s}\left([Q\circ^{0\mathbf{PT}_{\boldsymbol{\mathcal{A}}}}_{s} P]_{s}\right)
\tag{3}
\\&=
\mathrm{ip}^{(2,[1])\sharp}_{s}\left([Q]_{s}\right)
\circ^{0\mathbf{Pth}_{\boldsymbol{\mathcal{A}}^{(2)}}}_{s}
\mathrm{ip}^{(2,[1])\sharp}_{s}\left([P]_{s}\right)
\tag{4}
\\&=
\mathrm{ip}^{(2,X)@}_{s}\left(
\mathrm{CH}^{(2)}_{s}\left(
\mathrm{ip}^{(2,[1])\sharp}_{s}\left([Q]_{s}\right)
\right)\right)
\circ^{0\mathbf{Pth}_{\boldsymbol{\mathcal{A}}^{(2)}}}_{s}
\\&\qquad\qquad\qquad\qquad\qquad\qquad\qquad
\mathrm{ip}^{(2,X)@}_{s}\left(
\mathrm{CH}^{(2)}_{s}\left(
\mathrm{ip}^{(2,[1])\sharp}_{s}\left([P]_{s}\right)
\right)\right)
\tag{5}
\\&=
\mathrm{ip}^{(2,X)@}_{s}\left(
\mathrm{CH}^{(2)}_{s}\left(
\mathfrak{Q}^{(2)}
\right)\right)
\circ^{0\mathbf{Pth}_{\boldsymbol{\mathcal{A}}^{(2)}}}_{s}
\mathrm{ip}^{(2,X)@}_{s}\left(
\mathrm{CH}^{(2)}_{s}\left(
\mathfrak{P}^{(2)}
\right)\right).
\tag{6}
\end{align*}
\end{flushleft}

In the just stated chain of equalities, the first equality unravels the description of the second-order paths $\mathfrak{Q}^{(2)}$ and $\mathfrak{P}^{(2)}$; the second equality follows from Proposition~\ref{PDUIp}; the third equality follows from Proposition~\ref{PDIpDU}; the fourth equality follows from Proposition~\ref{PDUIp}; the fifth equality follows from Proposition~\ref{PDUIp}; finally, the last equality recovers the description of the second-order paths $\mathfrak{Q}^{(2)}$ and $\mathfrak{P}^{(2)}$.

This completes Case~(1).

Now, consider Case~(2), i.e., the case in which either $\mathfrak{Q}^{(2)}$ or $\mathfrak{P}^{(2)}$ are not $(2,[1])$-identity second-order paths.
In this case, the following chain of equalities holds
\begin{flushleft}
$\mathrm{ip}^{(2,X)@}_{s}\left(
\mathrm{CH}^{(2)}_{s}\left(
\mathfrak{Q}^{(2)}
\circ^{0\mathbf{Pth}_{\boldsymbol{\mathcal{A}}^{(2)}}}_{s}
\mathfrak{P}^{(2)}
\right)\right)$
\allowdisplaybreaks
\begin{align*}
&=
\mathrm{ip}^{(2,X)@}_{s}\left(
\mathrm{CH}^{(2)}_{s}\left(
\mathfrak{Q}^{(2)}\right)
\circ^{0\mathbf{T}_{\Sigma^{\boldsymbol{\mathcal{A}}^{(2)}}}(X)}_{s}
\mathrm{CH}^{(2)}_{s}\left(
\mathfrak{P}^{(2)}
\right)
\right)
\tag{1}
\\&=
\mathrm{ip}^{(2,X)@}_{s}\left(
\mathrm{CH}^{(2)}_{s}\left(
\mathfrak{Q}^{(2)}\right)
\right)
\circ^{0\mathbf{Pth}_{\boldsymbol{\mathcal{A}}^{(2)}}}_{s}
\mathrm{ip}^{(2,X)@}_{s}\left(
\mathrm{CH}^{(2)}_{s}\left(
\mathfrak{P}^{(2)}
\right)
\right).
\tag{2}
\end{align*}
\end{flushleft}

In the just stated chain of equalities, the first equality unravels the definition of the second-order Curry-Howard mapping, according to Definition~\ref{DDCH}. In this regard let us recall that, under the assumption that either $\mathfrak{P}^{(2)}$ or $\mathfrak{Q}^{(2)}$ is not a $(2,[1])$-identity second-order path, then $\mathfrak{Q}^{(2)}\circ^{0\mathbf{Pth}_{\boldsymbol{\mathcal{A}}^{(2)}}}_{s}
\mathfrak{P}^{(2)}$ is a coherent head-constant echelonless second-order path associated to the operation symbol of $0$-composition; the second equality follows from the fact that $\mathrm{ip}^{(2,X)@}$ is a $\Sigma^{\boldsymbol{\mathcal{A}}^{(2)}}$-homomorphism. In this regard, let us recall that the $0$-composition is well-defined in virtue of Proposition~\ref{PDIpDCH} and Corollary~\ref{CDCH}.

This completes Case~(2).

This completes the proof.
\end{proof}

We next prove that the normalized form of a second-order rewrite rule is the second-order rewrite rule itself.

\begin{restatable}{lemma}{LDIpDCHDEch}
\label{LDIpDCHDEch} Let $s$ be a sort in $S$ and let $\mathfrak{p}^{(2)}$ be a second-order rewrite rule in $\mathcal{A}^{(2)}_{s}$ then the following equality holds
\[
\mathrm{ip}^{(2,X)@}_{s}\left(
\mathrm{CH}^{(2)}_{s}\left(
\mathfrak{p}^{(2)\mathbf{Pth}_{\boldsymbol{\mathcal{A}}^{(2)}}}
\right)\right)
=
\mathfrak{p}^{(2)\mathbf{Pth}_{\boldsymbol{\mathcal{A}}^{(2)}}}.
\]
\end{restatable}
\begin{proof}
The following chain of equalities holds.
\allowdisplaybreaks
\begin{align*}
\mathrm{ip}^{(2,X)@}_{s}\left(
\mathrm{CH}^{(2)}_{s}\left(
\mathfrak{p}^{(2)\mathbf{Pth}_{\boldsymbol{\mathcal{A}}^{(2)}}}
\right)\right)&=
\mathrm{ip}^{(2,X)@}_{s}\left(
\mathrm{CH}^{(2)}_{s}\left(
\mathrm{ech}^{(2,\mathcal{A}^{(2)})}_{s}\left(
\mathfrak{p}^{(2)}
\right)
\right)\right)
\tag{1}
\\&=
\mathrm{ip}^{(2,X)@}_{s}\left(
\eta^{(2,\mathcal{A}^{(2)})}_{s}\left(
\mathfrak{p}^{(2)}
\right)\right)
\tag{2}
\\&=
\mathrm{ech}^{(2,\mathcal{A}^{(2)})}_{s}\left(
\mathfrak{p}^{(2)}
\right)
\tag{3}
\\&=
\mathfrak{p}^{(2)\mathbf{Pth}_{\boldsymbol{\mathcal{A}}^{(2)}}}.
\tag{4}
\end{align*}

In the just stated chain of equalities, the first equality unravels the interpretation of the constant operation symbol $\mathfrak{p}^{(2)}$ in the many-sorted partial $\Sigma^{\boldsymbol{\mathcal{A}}^{(2)}}$-algebra $\mathbf{Pth}_{\boldsymbol{\mathcal{A}}^{(2)}}$, according to Proposition~\ref{PDPthDCatAlg}; the second equality follows from Proposition~\ref{PDCHDA}; the third equality follows from Proposition~\ref{PDIpEch}; finally, the last equality recovers the interpretation of the constant operation symbol $\mathfrak{p}^{(2)}$ in the many-sorted partial $\Sigma^{\boldsymbol{\mathcal{A}}^{(2)}}$-algebra $\mathbf{Pth}_{\boldsymbol{\mathcal{A}}^{(2)}}$, according to Proposition~\ref{PDPthDCatAlg}.

This completes the proof.
\end{proof}

Finally we prove that the normalized form of a $1$-source of a second-order path is the $1$-source of the respective normalized form. The same also holds for the $1$-target.

\begin{restatable}{lemma}{LDIpDCHScUTgU}
\label{LDIpDCHScUTgU} Let $s$ be a sort in $S$ and let $\mathfrak{P}^{(2)}$ be a second-order path in $\mathrm{Pth}_{\boldsymbol{\mathcal{A}}^{(2)},s}$
then the following equality holds
\allowdisplaybreaks
\begin{align*}
\mathrm{ip}^{(2,X)@}_{s}\left(
\mathrm{CH}^{(2)}_{s}\left(
\mathrm{sc}^{1\mathbf{Pth}_{\boldsymbol{\mathcal{A}}^{(2)}}}_{s}\left(
\mathfrak{P}^{(2)}
\right)\right)\right)
&=
\mathrm{sc}^{1\mathbf{Pth}_{\boldsymbol{\mathcal{A}}^{(2)}}}_{s}\left(
\mathrm{ip}^{(2,X)@}_{s}\left(
\mathrm{CH}^{(2)}_{s}\left(
\mathfrak{P}^{(2)}
\right)\right)\right);
\\
\mathrm{ip}^{(2,X)@}_{s}\left(
\mathrm{CH}^{(2)}_{s}\left(
\mathrm{tg}^{1\mathbf{Pth}_{\boldsymbol{\mathcal{A}}^{(2)}}}_{s}\left(
\mathfrak{P}^{(2)}
\right)\right)\right)
&=
\mathrm{tg}^{1\mathbf{Pth}_{\boldsymbol{\mathcal{A}}^{(2)}}}_{s}\left(
\mathrm{ip}^{(2,X)@}_{s}\left(
\mathrm{CH}^{(2)}_{s}\left(
\mathfrak{P}^{(2)}
\right)\right)\right).
\end{align*}
\end{restatable}
\begin{proof}
We present the proof only for the first statement. The other one can be handled similarly.

The following chain of equalities holds
\begin{flushleft}
$\mathrm{ip}^{(2,X)@}_{s}\left(
\mathrm{CH}^{(2)}_{s}\left(
\mathrm{sc}^{1\mathbf{Pth}_{\boldsymbol{\mathcal{A}}^{(2)}}}_{s}\left(
\mathfrak{P}^{(2)}
\right)\right)\right)$
\allowdisplaybreaks
\begin{align*}
\qquad&=
\mathrm{ip}^{(2,X)@}_{s}\left(
\mathrm{CH}^{(2)}_{s}\left(
\mathrm{ip}^{(2,[1])\sharp}_{s}\left(
\mathrm{sc}^{([1],2)}_{s}\left(
\mathfrak{P}^{(2)}
\right)\right)\right)\right)
\tag{1}
\\&=
\mathrm{ip}^{(2,X)@}_{s}\left(
\eta^{(2,1)}_{s}\left(
\mathrm{CH}^{(1)\mathrm{m}}_{s}\left(
\mathrm{ip}^{([1],X)@}_{s}\left(
\mathrm{sc}^{([1],2)}_{s}\left(
\mathfrak{P}^{(2)}
\right)\right)\right)\right)\right)
\tag{2}
\\&=
\mathrm{ip}^{(2,[1])\sharp}_{s}\left(
\left[
\mathrm{CH}^{(1)\mathrm{m}}_{s}\left(
\mathrm{ip}^{([1],X)@}_{s}\left(
\mathrm{sc}^{([1],2)}_{s}\left(
\mathfrak{P}^{(2)}
\right)\right)\right)
\right]_{s}
\right)
\tag{3}
\\&=
\mathrm{ip}^{(2,[1])\sharp}_{s}\left(
\mathrm{sc}^{([1],2)}_{s}\left(
\mathfrak{P}^{(2)}
\right)\right)
\tag{4}
\\&=
\mathrm{ip}^{(2,[1])\sharp}_{s}\left(
\mathrm{sc}^{([1],2)}_{s}\left(
\mathrm{ip}^{(2,X)@}_{s}\left(
\mathrm{CH}^{(2)}_{s}\left(
\mathfrak{P}^{(2)}
\right)\right)
\right)\right)
\tag{5}
\\&=
\mathrm{sc}^{1\mathbf{Pth}_{\boldsymbol{\mathcal{A}}^{(2)}}}_{s}\left(
\mathrm{ip}^{(2,X)@}_{s}\left(
\mathrm{CH}^{(2)}_{s}\left(
\mathfrak{P}^{(2)}
\right)\right)
\right).
\tag{6}
\end{align*}
\end{flushleft}

In the just stated chain of equalities, the first equality unravels the description of the operation symbol $\mathrm{sc}^{1}_{s}$ in the many-sorted partial $\Sigma^{\boldsymbol{\mathcal{A}}^{(2)}}$-algebra $\mathbf{Pth}_{\boldsymbol{\mathcal{A}}^{(2)}}$, according to Proposition~\ref{PDPthDCatAlg}; the second equality follows from Proposition~\ref{PDCHDUId}; the third equality follows from Proposition~\ref{PDIpDU}; the fourth equality follows from Lemma~\ref{LWCong};
the fifth equality follows from Proposition~\ref{PDIpDCH} and Lemma~\ref{LDCH}; finally, the last equality recovers the description of the operation symbol $\mathrm{sc}^{1}_{s}$ in the many-sorted partial $\Sigma^{\boldsymbol{\mathcal{A}}^{(2)}}$-algebra $\mathbf{Pth}_{\boldsymbol{\mathcal{A}}^{(2)}}$, according to Proposition~\ref{PDPthCatAlg}.

This completes the proof.
\end{proof}

The following corollary states that one-step second-order paths are always normalized. In case there is just one derivation, this has been done, by necessity, in the leftmost innermost possible position.

\begin{restatable}{corollary}{CDIpDCHOneStep}
\label{CDIpDCHOneStep}
Let $s$ be a sort in $S$, and $\mathfrak{P}^{(2)}$ a second-order path in $\mathrm{Pth}_{\boldsymbol{\mathcal{A}}^{(2)},s}$. If $\mathfrak{P}^{(2)}$ is a one-step second-order path, then $\mathrm{ip}^{(2,X)@}_{s}(\mathrm{CH}^{(2)}_{s}(\mathfrak{P}^{(2)}))=\mathfrak{P}^{(2)}$.
\end{restatable}
\begin{proof}
It follows from Propositions~\ref{PDIpDCH} and~\ref{PDCHOneStep}.
\end{proof}

\section{
\texorpdfstring
{Interaction with the congruence $\Upsilon^{[1]}$ on $\mathbf{Pth}_{\boldsymbol{\mathcal{A}}^{(2)}}$}
{Interaction with Upsilon}
}

In this section we study the relation of a second-order path to its normalised version with respect to the $\Sigma^{\boldsymbol{\mathcal{A}}^{(2)}}$-congruence $\Upsilon^{[1]}$ on $\mathbf{Pth}_{\boldsymbol{\mathcal{A}}^{(2)}}$ introduced in Definition~\ref{DDUpsCong}. 

The following proposition proves that, for every sort $s\in S$, every $(2,[1])$-identity second-order path $\mathfrak{P}^{(2)}$ in $\mathrm{Pth}_{\boldsymbol{\mathcal{A}}^{(2)},s}$  is $\Upsilon^{[1]}_{s}$-related to its normalised second-order path.

\begin{restatable}{proposition}{PDIpUpsId}
\label{PDIpUpsId} Let $s$ be a sort in $S$ and let $\mathfrak{P}^{(2)}$ be a $(2,[1])$-identity second-order path in $\mathrm{Pth}_{\boldsymbol{\mathcal{A}}^{(2)},s}$, then
\[
\left(
\mathfrak{P}^{(2)},
\mathrm{ip}^{(2,X)@}_{s}\left(
\mathrm{CH}^{(2)}_{s}\left(
\mathfrak{P}^{(2)}
\right)
\right)
\right)
\in\Upsilon^{[1]}_{s}.
\]
\end{restatable}
\begin{proof}
It follows directly from Proposition~\ref{PDIpDU}.
\end{proof}

Next we show that if two second-order paths are related for the $\Sigma^{\boldsymbol{\mathcal{A}}^{(2)}}$-congruence $\Upsilon^{[1]}$, then so are their normalised versions.

\begin{restatable}{proposition}{PDIpUps}
\label{PDIpUps} Let $s$ be a sort in $S$ and $\mathfrak{P}^{(2)}, \mathfrak{Q}^{(2)}$ be second-order paths in $\mathrm{Pth}_{\boldsymbol{\mathcal{A}}^{(2)},s}$ satisfying that $(\mathfrak{P}^{(2)},\mathfrak{Q}^{(2)})\in \Upsilon^{(1)}_{s}$, then the following property holds
\[
\left(
\mathrm{ip}^{(2,X)@}_{s}\left(
\mathrm{CH}^{(2)}_{s}\left(
\mathfrak{P}^{(2)}
\right)\right)
,
\mathrm{ip}^{(2,X)@}_{s}\left(
\mathrm{CH}^{(2)}_{s}\left(
\mathfrak{Q}^{(2)}
\right)\right)
\right)
\in\Upsilon^{(1)}_{s}.
\]
\end{restatable}
\begin{proof}
We recall the definition of the relation $\Upsilon^{(1)}$ on $\mathrm{Pth}_{\boldsymbol{\mathcal{A}}^{(2)}}$, introduced in Definition~\ref{DDUps}. We will consider each case individually.

\begin{claim}\label{CDIpUpsI} Let $s$ be a sort in $S$ and $\mathfrak{P}^{(2)}$ a second-order path in $\mathrm{Pth}_{\boldsymbol{\mathcal{A}}^{(2)},s}$, then
\allowdisplaybreaks
\begin{multline*}
\left(
\mathrm{ip}^{(2,X)@}_{s}\left(
\mathrm{CH}^{(2)}_{s}\left(
\mathfrak{P}^{(2)}
\right)
\right),
\right.
\\
\left.
\mathrm{ip}^{(2,X)@}_{s}\left(
\mathrm{CH}^{(2)}_{s}\left(
\mathfrak{P}^{(2)}\circ^{0\mathbf{Pth}_{\boldsymbol{\mathcal{A}}^{(2)}}}_{s}
\mathrm{sc}^{0\mathbf{Pth}_{\boldsymbol{\mathcal{A}}^{(2)}}}_{s}\left(
\mathfrak{P}^{(2)}
\right)
\right)
\right)
\right)\in \Upsilon^{[1]}_{s}.
\tag{1}
\end{multline*}
\end{claim}

The following chain of equalities holds
\begin{flushleft}
$\mathrm{ip}^{(2,X)@}_{s}\left(
\mathrm{CH}^{(2)}_{s}\left(
\mathfrak{P}^{(2)}\circ^{0\mathbf{Pth}_{\boldsymbol{\mathcal{A}}^{(2)}}}_{s}
\mathrm{sc}^{0\mathbf{Pth}_{\boldsymbol{\mathcal{A}}^{(2)}}}_{s}\left(
\mathfrak{P}^{(2)}
\right)
\right)
\right)
$
\allowdisplaybreaks
\begin{align*}
&=
\mathrm{ip}^{(2,X)@}_{s}\left(
\mathrm{CH}^{(2)}_{s}\left(
\mathfrak{P}^{(2)}\right)
\right)
\circ^{0\mathbf{Pth}_{\boldsymbol{\mathcal{A}}^{(2)}}}_{s}
\mathrm{ip}^{(2,X)@}_{s}\left(
\mathrm{CH}^{(2)}_{s}\left(
\mathrm{sc}^{0\mathbf{Pth}_{\boldsymbol{\mathcal{A}}^{(2)}}}_{s}\left(
\mathfrak{P}^{(2)}
\right)
\right)
\right)
\tag{1}
\\&=
\mathrm{ip}^{(2,X)@}_{s}\left(
\mathrm{CH}^{(2)}_{s}\left(
\mathfrak{P}^{(2)}\right)
\right)
\circ^{0\mathbf{Pth}_{\boldsymbol{\mathcal{A}}^{(2)}}}_{s}
\mathrm{sc}^{0\mathbf{Pth}_{\boldsymbol{\mathcal{A}}^{(2)}}}_{s}\left(
\mathfrak{P}^{(2)}
\right)
\tag{2}
\\&=
\mathrm{ip}^{(2,X)@}_{s}\left(
\mathrm{CH}^{(2)}_{s}\left(
\mathfrak{P}^{(2)}\right)
\right)
\circ^{0\mathbf{Pth}_{\boldsymbol{\mathcal{A}}^{(2)}}}_{s}
\mathrm{ip}^{(2,0)\sharp}_{s}\left(
\mathrm{sc}^{(0,2)}_{s}\left(
\mathfrak{P}^{(2)}
\right)\right)
\tag{3}
\\&=
\mathrm{ip}^{(2,X)@}_{s}\left(
\mathrm{CH}^{(2)}_{s}\left(
\mathfrak{P}^{(2)}\right)
\right)
\circ^{0\mathbf{Pth}_{\boldsymbol{\mathcal{A}}^{(2)}}}_{s}
\\&\qquad\qquad\qquad\qquad\qquad\qquad\qquad
\mathrm{ip}^{(2,0)\sharp}_{s}\left(
\mathrm{sc}^{(0,2)}_{s}\left(
\mathrm{ip}^{(2,X)@}_{s}\left(
\mathrm{CH}^{(2)}_{s}\left(
\mathfrak{P}^{(2)}\right)
\right)
\right)\right)
\tag{4}
\\&=
\mathrm{ip}^{(2,X)@}_{s}\left(
\mathrm{CH}^{(2)}_{s}\left(
\mathfrak{P}^{(2)}\right)
\right)
\circ^{0\mathbf{Pth}_{\boldsymbol{\mathcal{A}}^{(2)}}}_{s}
\mathrm{sc}^{0\mathbf{Pth}_{\boldsymbol{\mathcal{A}}^{(2)}}}_{s}\left(
\mathrm{ip}^{(2,X)@}_{s}\left(
\mathrm{CH}^{(2)}_{s}\left(
\mathfrak{P}^{(2)}\right)
\right)
\right).
\tag{5}
\end{align*}
\end{flushleft}

In the just stated chain of equalities, the first equality follows from Lemma~\ref{LDIpDCHCompZ}; the second equality follows from Claim~\ref{CDPthCatAlgScZ} and Proposition~\ref{PDIpDU}; the third equality unravels the description of the $0$-source operation in $\mathbf{Pth}_{\boldsymbol{\mathcal{A}}^{(2)}}$, according to Proposition~\ref{PDPthCatAlg}; the fourth equality follows from Proposition~\ref{PDIpDCH} and Corollary~\ref{CDCH}; finally, the last equality recovers the description of the $0$-source operation in $\mathbf{Pth}_{\boldsymbol{\mathcal{A}}^{(2)}}$, according to Proposition~\ref{PDPthCatAlg}.

Thus, the claim for this case follows because the pair under consideration falls under Case~(i), according to Definition~\ref{DDUps}.

This finishes the proof of Claim~\ref{CDIpUpsI}.

\begin{claim}\label{CDIpUpsII} Let $s$ be a sort in $S$ and $\mathfrak{P}^{(2)}$ a second-order path in $\mathrm{Pth}_{\boldsymbol{\mathcal{A}}^{(2)},s}$, then
\allowdisplaybreaks
\begin{multline*}
\left(
\mathrm{ip}^{(2,X)@}_{s}\left(
\mathrm{CH}^{(2)}_{s}\left(
\mathfrak{P}^{(2)}
\right)
\right),
\right.
\\
\left.
\mathrm{ip}^{(2,X)@}_{s}\left(
\mathrm{CH}^{(2)}_{s}\left(
\mathrm{tg}^{0\mathbf{Pth}_{\boldsymbol{\mathcal{A}}^{(2)}}}_{s}\left(
\mathfrak{P}^{(2)}
\right)
\circ^{0\mathbf{Pth}_{\boldsymbol{\mathcal{A}}^{(2)}}}_{s}
\mathfrak{P}^{(2)}
\right)
\right)
\right)\in \Upsilon^{[1]}_{s}.
\tag{2}
\end{multline*}
\end{claim}

This claim follows by a similar argument to that used for the proof of Claim~\ref{CDIpUpsI}.

\begin{claim}\label{CDIpUpsIII} Let $s$ be a sort in $S$ and $\mathfrak{P}^{(2)}, \mathfrak{Q}^{(2)}$ and $\mathfrak{R}^{(2)}$ be second-order paths in $\mathrm{Pth}_{\boldsymbol{\mathcal{A}}^{(2)},s}$ with
\begin{align*}
\mathrm{sc}^{(0,2)}_{s}\left(
\mathfrak{R}^{(2)}
\right)
&=
\mathrm{tg}^{(0,2)}_{s}\left(
\mathfrak{R}^{(2)}
\right);
&
\mathrm{sc}^{(0,2)}_{s}\left(
\mathfrak{Q}^{(2)}
\right)
&=
\mathrm{tg}^{(0,2)}_{s}\left(
\mathfrak{P}^{(2)}
\right).
\end{align*}
then the following property holds
\begin{multline*}
\left(
\mathrm{ip}^{(2,X)@}_{s}\left(
\mathrm{CH}^{(2)}_{s}\left(
\mathfrak{R}^{(2)}
\circ^{0\mathbf{Pth}_{\boldsymbol{\mathcal{A}}^{(2)}}}_{s}
\left(
\mathfrak{Q}^{(2)}
\circ^{0\mathbf{Pth}_{\boldsymbol{\mathcal{A}}^{(2)}}}_{s}
\mathfrak{P}^{(2)}
\right)
\right)
\right),
\right.
\\
\left.
\mathrm{ip}^{(2,X)@}_{s}\left(
\mathrm{CH}^{(2)}_{s}\left(
\left(
\mathfrak{R}^{(2)}
\circ^{0\mathbf{Pth}_{\boldsymbol{\mathcal{A}}^{(2)}}}_{s}
\mathfrak{Q}^{(2)}
\right)
\circ^{0\mathbf{Pth}_{\boldsymbol{\mathcal{A}}^{(2)}}}_{s}
\mathfrak{P}^{(2)}
\right)
\right)
\right)\in \Upsilon^{[1]}_{s}.
\tag{3}
\end{multline*}
\end{claim}

Consider the left hand side of the previous pair. In virtue of Lemma~\ref{LDIpDCHCompZ}, the following equality holds
\allowdisplaybreaks
\begin{multline*}
\mathrm{ip}^{(2,X)@}_{s}\left(
\mathrm{CH}^{(2)}_{s}\left(
\mathfrak{R}^{(2)}
\circ^{0\mathbf{Pth}_{\boldsymbol{\mathcal{A}}^{(2)}}}_{s}
\left(
\mathfrak{Q}^{(2)}
\circ^{0\mathbf{Pth}_{\boldsymbol{\mathcal{A}}^{(2)}}}_{s}
\mathfrak{P}^{(2)}
\right)
\right)
\right)
\\=
\mathrm{ip}^{(2,X)@}_{s}\left(
\mathrm{CH}^{(2)}_{s}\left(
\mathfrak{R}^{(2)}
\right)\right)
\circ^{0\mathbf{Pth}_{\boldsymbol{\mathcal{A}}^{(2)}}}_{s}
\left(
\mathrm{ip}^{(2,X)@}_{s}\left(
\mathrm{CH}^{(2)}_{s}\left(
\mathfrak{Q}^{(2)}
\right)\right)
\right.
\\
\left.
\circ^{0\mathbf{Pth}_{\boldsymbol{\mathcal{A}}^{(2)}}}_{s}
\mathrm{ip}^{(2,X)@}_{s}\left(
\mathrm{CH}^{(2)}_{s}\left(
\mathfrak{P}^{(2)}
\right)\right)
\right).
\end{multline*}

Now, consider the right hand side of the previous pair. In virtue of Lemma~\ref{LDIpDCHCompZ}, the following equality holds
\allowdisplaybreaks
\begin{multline*}
\mathrm{ip}^{(2,X)@}_{s}\left(
\mathrm{CH}^{(2)}_{s}\left(
\left(
\mathfrak{R}^{(2)}
\circ^{0\mathbf{Pth}_{\boldsymbol{\mathcal{A}}^{(2)}}}_{s}
\mathfrak{Q}^{(2)}
\right)
\circ^{0\mathbf{Pth}_{\boldsymbol{\mathcal{A}}^{(2)}}}_{s}
\mathfrak{P}^{(2)}
\right)
\right)
\\=
\left(
\mathrm{ip}^{(2,X)@}_{s}\left(
\mathrm{CH}^{(2)}_{s}\left(
\mathfrak{R}^{(2)}
\right)\right)
\circ^{0\mathbf{Pth}_{\boldsymbol{\mathcal{A}}^{(2)}}}_{s}
\mathrm{ip}^{(2,X)@}_{s}\left(
\mathrm{CH}^{(2)}_{s}\left(
\mathfrak{Q}^{(2)}
\right)
\right)
\right)
\\
\circ^{0\mathbf{Pth}_{\boldsymbol{\mathcal{A}}^{(2)}}}_{s}
\mathrm{ip}^{(2,X)@}_{s}\left(
\mathrm{CH}^{(2)}_{s}\left(
\mathfrak{P}^{(2)}
\right)\right).
\end{multline*}

Thus, the claim for this case follows because the pair under consideration falls under Case~(iii), according to Definition~\ref{DDUps}.

This finishes the proof of Claim~\ref{CDIpUpsIII}.

\begin{claim}\label{CDIpUpsIV} Let $(\mathbf{s},s)$ be a pair in $S^{\star}\times S$, $\sigma\in \Sigma_{\mathbf{s},s,}$ and $(\mathfrak{P}^{(2)}_{j})_{j\in\bb{\mathbf{s}}}, (\mathfrak{P}^{(2)}_{j})_{j\in\bb{\mathbf{s}}}$ two families of second-order paths in $\mathrm{Pth}_{\boldsymbol{\mathcal{A}}^{(2)},\mathbf{s}}$ satisfying that, for every $j\in\bb{\mathbf{s}}$, 
\[
\mathrm{sc}^{(0,2)}_{s_{j}}\left(
\mathfrak{Q}^{(2)}_{j}
\right)
=
\mathrm{tg}^{(0,2)}_{s_{j}}\left(
\mathfrak{P}^{(2)}_{j}
\right)
\]
then the following property holds
\begin{multline*}
\left(
\mathrm{ip}^{(2,X)@}_{s}\left(
\mathrm{CH}^{(2)}_{s}\left(
\sigma^{\mathbf{Pth}_{\boldsymbol{\mathcal{A}}^{(2)}}}\left(
\left(
\mathfrak{Q}^{(2)}_{j}
\circ^{0\mathbf{Pth}_{\boldsymbol{\mathcal{A}}^{(2)}}}_{s_{j}}
\mathfrak{P}^{(2)}_{j}
\right)_{j\in\bb{\mathbf{s}}}
\right)
\right)
\right),
\right.
\\
\mathrm{ip}^{(2,X)@}_{s}\left(
\mathrm{CH}^{(2)}_{s}\left(
\sigma^{\mathbf{Pth}_{\boldsymbol{\mathcal{A}}^{(2)}}}
\left(
\left(
\mathfrak{Q}^{(2)}_{j}
\right)_{j\in\bb{\mathbf{s}}}
\right)
\right.\right.
\\
\left.\left.\left.
\circ^{0\mathbf{Pth}_{\boldsymbol{\mathcal{A}}^{(2)}}}_{s}
\sigma^{\mathbf{Pth}_{\boldsymbol{\mathcal{A}}^{(2)}}}
\left(
\left(
\mathfrak{P}^{(2)}_{j}
\right)_{j\in\bb{\mathbf{s}}}
\right)
\right)
\right)
\right)\in \Upsilon^{[1]}_{s}.
\tag{4}
\end{multline*}
\end{claim}

Consider the left hand side of the previous pair. The following chain of equalities holds
\begin{flushleft}
$\mathrm{ip}^{(2,X)@}_{s}\left(
\mathrm{CH}^{(2)}_{s}\left(
\sigma^{\mathbf{Pth}_{\boldsymbol{\mathcal{A}}^{(2)}}}\left(
\left(
\mathfrak{Q}^{(2)}_{j}
\circ^{0\mathbf{Pth}_{\boldsymbol{\mathcal{A}}^{(2)}}}_{s_{j}}
\mathfrak{P}^{(2)}_{j}
\right)_{j\in\bb{\mathbf{s}}}
\right)
\right)
\right)$
\allowdisplaybreaks
\begin{align*}
&=
\sigma^{\mathbf{Pth}_{\boldsymbol{\mathcal{A}}^{(2)}}}\left(
\left(
\mathrm{ip}^{(2,X)@}_{s_{j}}\left(
\mathrm{CH}^{(2)}_{s_{j}}\left(
\mathfrak{Q}^{(2)}_{j}
\circ^{0\mathbf{Pth}_{\boldsymbol{\mathcal{A}}^{(2)}}}_{s_{j}}
\mathfrak{P}^{(2)}_{j}
\right)
\right)
\right)_{j\in\bb{\mathbf{s}}}
\right)
\tag{1}
\\&=
\sigma^{\mathbf{Pth}_{\boldsymbol{\mathcal{A}}^{(2)}}}\left(
\left(
\mathrm{ip}^{(2,X)@}_{s_{j}}\left(
\mathrm{CH}^{(2)}_{s_{j}}\left(
\mathfrak{Q}^{(2)}_{j}
\right)\right)
\circ^{0\mathbf{Pth}_{\boldsymbol{\mathcal{A}}^{(2)}}}_{s_{j}}
\right.\right.
\\&\qquad\qquad\qquad\qquad\qquad\qquad\qquad
\left.\left.
\mathrm{ip}^{(2,X)@}_{s_{j}}\left(
\mathrm{CH}^{(2)}_{s_{j}}\left(
\mathfrak{P}^{(2)}_{j}
\right)
\right)
\right)_{j\in\bb{\mathbf{s}}}
\right).
\tag{2}
\end{align*}
\end{flushleft}

In the just stated chain of equalities, the first equality follows from the fact that $\mathrm{CH}^{(2)}$ and $\mathrm{ip}^{(2,X)@}$ are $\Sigma$-homomorphism in virtue of Proposition~\ref{PDCHHom} and Definition~\ref{DDIp}, respectively; finally the last equality follows from Lemma~\ref{LDIpDCHCompZ}.

Now, consider the right hand side of the pair in Equation (4) above. The following chain of equalities holds
\begin{flushleft}
$\mathrm{ip}^{(2,X)@}_{s}\left(
\mathrm{CH}^{(2)}_{s}\left(
\sigma^{\mathbf{Pth}_{\boldsymbol{\mathcal{A}}^{(2)}}}
\left(
\left(
\mathfrak{Q}^{(2)}_{j}
\right)_{j\in\bb{\mathbf{s}}}
\right)
\circ^{0\mathbf{Pth}_{\boldsymbol{\mathcal{A}}^{(2)}}}_{s}
\sigma^{\mathbf{Pth}_{\boldsymbol{\mathcal{A}}^{(2)}}}
\left(
\left(
\mathfrak{P}^{(2)}_{j}
\right)_{j\in\bb{\mathbf{s}}}
\right)
\right)
\right)$
\allowdisplaybreaks
\begin{align*}
\qquad&=
\mathrm{ip}^{(2,X)@}_{s}\left(
\mathrm{CH}^{(2)}_{s}\left(
\sigma^{\mathbf{Pth}_{\boldsymbol{\mathcal{A}}^{(2)}}}
\left(
\left(
\mathfrak{Q}^{(2)}_{j}
\right)_{j\in\bb{\mathbf{s}}}
\right)
\right)\right)
\circ^{0\mathbf{Pth}_{\boldsymbol{\mathcal{A}}^{(2)}}}_{s}
\\&\qquad\qquad\qquad\qquad\qquad\qquad
\mathrm{ip}^{(2,X)@}_{s}\left(
\mathrm{CH}^{(2)}_{s}\left(
\sigma^{\mathbf{Pth}_{\boldsymbol{\mathcal{A}}^{(2)}}}
\left(
\left(
\mathfrak{P}^{(2)}_{j}
\right)_{j\in\bb{\mathbf{s}}}
\right)
\right)
\right)
\tag{1}
\\&=
\sigma^{\mathbf{Pth}_{\boldsymbol{\mathcal{A}}^{(2)}}}
\left(
\left(
\mathrm{ip}^{(2,X)@}_{s_{j}}\left(
\mathrm{CH}^{(2)}_{s_{j}}\left(
\mathfrak{Q}^{(2)}_{j}
\right)\right)
\right)_{j\in\bb{\mathbf{s}}}
\right)
\circ^{0\mathbf{Pth}_{\boldsymbol{\mathcal{A}}^{(2)}}}_{s}
\\&\qquad\qquad\qquad\qquad\qquad\qquad
\sigma^{\mathbf{Pth}_{\boldsymbol{\mathcal{A}}^{(2)}}}
\left(
\left(
\mathrm{ip}^{(2,X)@}_{s_{j}}\left(
\mathrm{CH}^{(2)}_{s_{j}}\left(
\mathfrak{P}^{(2)}_{j}
\right)\right)
\right)_{j\in\bb{\mathbf{s}}}
\right).
\tag{2}
\end{align*}
\end{flushleft}

In the just stated chain of equalities, the first equality follows from Lemma~\ref{LDIpDCHCompZ}; finally the last equality follows  from the fact that $\mathrm{CH}^{(2)}$ and $\mathrm{ip}^{(2,X)@}$ are $\Sigma$-homomorphisms in virtue of Proposition~\ref{PDCHHom} and Definition~\ref{DDIp}, respectively.

Thus, the claim for this case follows because the pair under consideration falls under Case~(iv), according to Definition~\ref{DDUps}.

This finishes the proof of Claim~\ref{CDIpUpsIV}.

This completes the proof.
\end{proof}

\section{Order issues}

In this section we study how the $\mathrm{ip}^{(2,X)@}$ mapping relates with the order on terms. Next proposition states that if the image under the second-order Curry-Howard mapping of a second-order path is smaller, with respect to the subterm order than the image under the second-order Curry-Howard mapping of another second-order path, then the respective normalised second-order paths satisfy the same inequality in the order for second-order paths introduced in Definition~\ref{DDOrd}.

\begin{restatable}{proposition}{PDIpDCHOrd}
\label{PDIpDCHOrd} Let $s$ and $t$ be sorts in $S$ and let $\mathfrak{Q}^{(2)}$ and $\mathfrak{P}^{(2)}$ be second-order paths in $\mathrm{Pth}_{\boldsymbol{\mathcal{A}}^{(2)},t}$ and $\mathrm{Pth}_{\boldsymbol{\mathcal{A}}^{(2)},s}$, respectively, satisfying that 
\[
\left(\mathrm{CH}^{(2)}_{t}\left(\mathfrak{Q}^{(2)}\right),t\right)
\leq_{\mathbf{T}_{\Sigma^{\boldsymbol{\mathcal{A}}^{(2)}}}(X)}
\left(\mathrm{CH}^{(2)}_{s}\left(\mathfrak{P}^{(2)}\right),s\right).
\]
Then the following inequality holds
\[
\left(
\mathrm{ip}^{(2,X)@}_{t}\left(
\mathrm{CH}^{(2)}_{t}\left(
\mathfrak{Q}^{(2)}
\right)\right),t\right)
\leq_{\mathbf{Pth}_{\boldsymbol{\mathcal{A}}^{(2)}}}
\left(
\mathrm{ip}^{(2,X)@}_{s}\left(
\mathrm{CH}^{(2)}_{s}\left(
\mathfrak{P}^{(2)}
\right)\right),s\right).
\]
\end{restatable}
\begin{proof}
We consider, first of all, some extreme cases that can be easily handled and will help us to take further assumptions.

We begin with the case of $s=t$ and $\mathrm{CH}^{(2)}_{t}(\mathfrak{Q}^{(2)})=\mathrm{CH}^{(2)}_{s}(\mathfrak{P}^{(2)})$. In this case, by reflexivity of $\leq_{\mathbf{Pth}_{\boldsymbol{\mathcal{A}}^{(2)}}}$ we have that
\[
\left(
\mathrm{ip}^{(2,X)@}_{t}\left(
\mathrm{CH}^{(2)}_{t}\left(
\mathfrak{Q}^{(2)}
\right)\right),t\right)
\leq_{\mathbf{Pth}_{\boldsymbol{\mathcal{A}}^{(2)}}}
\left(
\mathrm{ip}^{(2,X)@}_{s}\left(
\mathrm{CH}^{(2)}_{s}\left(
\mathfrak{P}^{(2)}
\right)\right),s\right).
\]

Therefore we can consider the following extra assumption
\begin{assumption}\label{ADIpDCHOrdI} We can assume that 
\[
\left(\mathrm{CH}^{(2)}_{t}\left(\mathfrak{Q}^{(2)}\right),t\right)
<_{\mathbf{T}_{\Sigma^{\boldsymbol{\mathcal{A}}^{(2)}}}(X)}
\left(\mathrm{CH}^{(2)}_{s}\left(\mathfrak{P}^{(2)}\right),s\right).
\]
\end{assumption}

Now, we consider the case in which $\mathfrak{P}^{(2)}$ is a $(2,[1])$-identity second-order path. Then for some path term $P$ in $\mathrm{PT}_{\boldsymbol{\mathcal{A}}, s}$ it is the case that 
\[
\mathfrak{P}^{(2)}=\mathrm{ip}^{(2,[1])\sharp}_{s}\left(
\left[
P
\right]_{s}
\right).
\]

In this case, according to Proposition~\ref{PDCHDUId}, we have that 
\[
\mathrm{CH}^{(2)}_{s}\left(
\mathfrak{P}^{(2)}
\right)
=
\eta^{(2,1)\sharp}_{s}\left(
\mathrm{CH}^{(1)}_{s}\left(
\mathrm{ip}^{(1,X)@}_{s}\left(
P
\right)\right)
\right).
\]

Note that, from Proposition~\ref{PDCHDUId}, also follows that  $\mathrm{CH}^{(2)}_{s}\left(
\mathfrak{P}^{(2)}
\right)$ is a term in $\eta^{(2,1)\sharp}[\mathrm{PT}_{\boldsymbol{\mathcal{A}}}]_{s}$. Since $\mathrm{CH}^{(2)}_{t}(\mathfrak{Q}^{(2)})$ is a subterm of $\mathrm{CH}^{(2)}_{s}(\mathfrak{P}^{(2)})$ of sort $t$, we conclude that $\mathrm{CH}^{(2)}_{t}(\mathfrak{Q}^{(2)})$ is also a term in $\eta^{(2,1)\sharp}[\mathrm{PT}_{\boldsymbol{\mathcal{A}}}]_{t}$. Again, by Proposition~\ref{PDCHDUId}, we conclude that $\mathfrak{Q}^{(2)}$ is a $(2,[1])$-identity second-order path. 

Let $Q$ be a path term in $\mathrm{PT}_{\boldsymbol{\mathcal{A}}, t}$ for which 
\[
\mathfrak{Q}^{(2)}=\mathrm{ip}^{(2,[1])\sharp}_{t}\left(
\left[
Q
\right]_{t}
\right).
\]

In this case, according to Proposition~\ref{PDCHDUId}, we have that 
\[
\mathrm{CH}^{(2)}_{t}\left(
\mathfrak{Q}^{(2)}
\right)
=
\eta^{(2,1)\sharp}_{t}\left(
\mathrm{CH}^{(1)}_{t}\left(
\mathrm{ip}^{(1,X)@}_{t}\left(
Q
\right)\right)
\right).
\]

Since, by assumption, we have that 
$
(\mathrm{CH}^{(2)}_{t}(\mathfrak{Q}^{(2)}),t)
\leq_{\mathbf{T}_{\Sigma^{\boldsymbol{\mathcal{A}}^{(2)}}}(X)}
(\mathrm{CH}^{(2)}_{s}(\mathfrak{P}^{(2)}),s)
$,  we conclude that 
\[
\left(
\mathrm{CH}^{(1)}_{t}\left(
\mathrm{ip}^{(1,X)@}_{t}\left(
Q
\right)\right),t\right)
\leq_{\mathbf{PT}_{\boldsymbol{\mathcal{A}}}
}
\left(
\mathrm{CH}^{(1)}_{s}\left(
\mathrm{ip}^{(1,X)@}_{s}\left(
P
\right)\right),s\right).
\]

It follows from Corollary~\ref{CPTQOrd} that 
\[
\left(
\left[Q\right]_{t},t
\right)
\leq_{[\mathbf{PT}_{\boldsymbol{\mathcal{A}}}]}
\left(
\left[P\right]_{s},s
\right).
\]

Following Definition~\ref{DDOrd}, we conclude that 
\[
\left(\mathfrak{Q}^{(2)},t\right)
\leq_{\mathbf{Pth}_{\boldsymbol{\mathcal{A}}^{(2)}}}
\left(\mathfrak{P}^{(2)},s\right).
\]

Finally, note that, according to Definition~\ref{DDIp} and Corollary~\ref{CDCHUId}, we have that
\begin{align*}
\mathrm{ip}^{(2,X)@}_{t}\left(
\mathrm{CH}_{t}\left(
\mathfrak{Q}^{(2)}
\right)\right)
&=\mathfrak{Q}^{(2)};
&
\mathrm{ip}^{(2,X)@}_{s}\left(
\mathrm{CH}_{s}\left(
\mathfrak{P}^{(2)}
\right)\right)
&=\mathfrak{P}^{(2)}.
\end{align*}

This completes the case of $\mathfrak{P}^{(2)}$ being a $(2,[1])$-identity second-order path.  It also allow us to take the following extra assumption.

\begin{assumption}\label{ADIpDCHOrdII} In what follows we can assume that $\mathfrak{P}^{(2)}$ is not a $(2,[1])$-identity second-order path.
\end{assumption}

We prove the general case by Artinian induction on $(\coprod\mathrm{Pth}_{\boldsymbol{\mathcal{A}}^{(2)}}, \leq_{\mathbf{Pth}_{\boldsymbol{\mathcal{A}}^{(2)}}})$ with respect to $(\mathfrak{P}^{(2)},s)$.  

\textsf{Base step of the induction.}

Assume that $(\mathfrak{P}^{(2)},s)$ is a minimal element in $(\coprod\mathrm{Pth}_{\boldsymbol{\mathcal{A}}^{(2)}}, \leq_{\mathbf{Pth}_{\boldsymbol{\mathcal{A}}^{(2)}}})$. Then according to Proposition~\ref{PDMinimal} and taking into account Assumption~\ref{ADIpDCHOrdII}, we are only left to prove the case of $\mathfrak{P}^{(2)}$ being  a  second-order echelon.

Now, if $\mathfrak{P}^{(2)}$ is a  second-order echelon associated to the second-order rewrite rule $\mathfrak{p}^{(2)}$ in $\mathcal{A}^{(2)}_{s}$, then following Definition~\ref{DDCH}, the image of the second-order Curry-Howard at $\mathfrak{P}^{(2)}$ is given by
\[
\mathrm{CH}^{(2)}_{s}\left(\mathfrak{P}^{(2)}\right)=\mathfrak{p}^{(2)\mathbf{T}_{\Sigma^{\boldsymbol{\mathcal{A}}^{(2)}}}(X)}.
\]

Note that $\mathfrak{p}^{(2)\mathbf{T}_{\Sigma^{\boldsymbol{\mathcal{A}}^{(2)}}}(X)}$ is a constant in $\mathbf{T}_{\Sigma^{\boldsymbol{\mathcal{A}}^{(2)}}}(X)_{s}$, i.e., a term in $\mathrm{B}^{0}_{\Sigma^{\boldsymbol{\mathcal{A}}^{(2)}}}(X)_{s}$. Therefore its only subterm is $\mathfrak{p}^{(2)\mathbf{T}_{\Sigma^{\boldsymbol{\mathcal{A}}^{(2)}}}(X)}$ itself. 

By assumption, we have that $
(\mathrm{CH}^{(2)}_{t}(\mathfrak{Q}^{(2)}),t)
\leq_{\mathbf{T}_{\Sigma^{\boldsymbol{\mathcal{A}}^{(2)}}}(X)}
(\mathrm{CH}^{(2)}_{s}(\mathfrak{P}^{(2)}),s)
$, therefore  we conclude that $t=s$  and $\mathrm{CH}^{(2)}_{t}(\mathfrak{Q}^{(2)})=\mathfrak{p}^{(2)\mathbf{T}_{\Sigma^{\boldsymbol{\mathcal{A}}^{(2)}}}(X)}$.

By reflexivity of $\leq_{\mathbf{Pth}_{\boldsymbol{\mathcal{A}}^{(2)}}}$ we conclude that 
\[
\left(
\mathrm{ip}^{(2,X)@}_{t}\left(
\mathrm{CH}^{(2)}_{t}\left(
\mathfrak{Q}^{(2)}
\right)\right),t\right)
\leq_{\mathbf{Pth}_{\boldsymbol{\mathcal{A}}^{(2)}}}
\left(
\mathrm{ip}^{(2,X)@}_{s}\left(
\mathrm{CH}^{(2)}_{s}\left(
\mathfrak{P}^{(2)}
\right)\right),s\right).
\]

This completes the case of $\mathfrak{P}^{(2)}$ being a  second-order echelon.

This complete the base step.

\textsf{Inductive step of the induction.}

Let $(\mathfrak{P}^{(2)},s)$ be a non-minimal element of $(\coprod\mathrm{Pth}_{\boldsymbol{\mathcal{A}}^{(2)}}, \leq_{\mathbf{Pth}_{\boldsymbol{\mathcal{A}}^{(2)}}})$. Let us suppose that, for every sort $u\in S$ and every second-order path $\mathfrak{R}^{(2)}$ in $\mathrm{Pth}_{\boldsymbol{\mathcal{A}}^{(2)},u}$, if $(\mathfrak{R}^{(2)},u)\leq_{\mathbf{Pth}_{\boldsymbol{\mathcal{A}}^{(2)}}}(\mathfrak{P}^{(2)},s)$ then the statement holds for $\mathfrak{R}^{(2)}$, i.e., for every sort $w\in S$, if $\mathfrak{S}^{(2)}$ is a second-order path in $\mathrm{Pth}_{\boldsymbol{\mathcal{A}}^{(2)},w}$ satisfying that 
\[
\left(
\mathrm{CH}^{(2)}_{w}\left(\mathfrak{S}^{(2)}\right),w
\right)
\leq_{\mathbf{T}_{\Sigma^{\boldsymbol{\mathcal{A}}^{(2)}}}(X)}
\left(
\mathrm{CH}^{(2)}_{u}\left(\mathfrak{R}^{(2)}\right),u
\right),
\]
then the following inequality holds
\[
\left(
\mathrm{ip}^{(2,X)@}_{w}\left(
\mathrm{CH}^{(2)}_{w}\left(
\mathfrak{S}^{(2)}
\right)\right),w\right)
\leq_{\mathbf{Pth}_{\boldsymbol{\mathcal{A}}^{(2)}}}
\left(
\mathrm{ip}^{(2,X)@}_{u}\left(
\mathrm{CH}^{(2)}_{u}\left(
\mathfrak{R}^{(2)}
\right)\right),u\right).
\]

Since $(\mathfrak{P}^{(2)},s)$ is a non-minimal element of $(\coprod\mathrm{Pth}_{\boldsymbol{\mathcal{A}}^{(2)}}, \leq_{\mathbf{Pth}_{\boldsymbol{\mathcal{A}}^{(2)}}})$ and taking into account Assumption~\ref{ADIpDCHOrdII}, we have, by Lemma~\ref{LDOrdI}, that $\mathfrak{P}^{(2)}$ is either (1) a second-order path of length strictly greater than one containing at least one  second-order echelon or (2) an echelonless secon-order path.

If~(1), then let $i\in\bb{\mathfrak{P}^{(2)}}$ be the first index for which the one-step subpath $\mathfrak{P}^{(2),i,i}$ is a  second-order echelon. We distinguish two cases accordingly (1.1) $i=0$ and (1.2) $i>0$.

If~(1.1) $i=0$, i.e., if $\mathfrak{P}^{(2)}$ has its first  second-order echelon on its first step, then according to Definition~\ref{DDCH}, we have that 
\[
\mathrm{CH}^{(2)}_{s}\left(\mathfrak{P}^{(2)}\right)
=
\mathrm{CH}^{(2)}_{s}\left(\mathfrak{P}^{(2),1,\bb{\mathfrak{P}^{(2)}}-1}\right)
\circ^{1\mathbf{T}_{\Sigma^{\boldsymbol{\mathcal{A}}^{(2)}}}(X)}_{s}
\mathrm{CH}^{(2)}_{s}\left(\mathfrak{P}^{(2),0,0}\right).
\]

Taking into account Proposition~\ref{PDIpDCH}, we have that $\mathrm{ip}^{(2,X)@}_{s}(\mathrm{CH}^{(2)}_{s}(\mathfrak{P}^{(2)}))$ is a second-order path in $[\mathfrak{P}^{(2)}]_{s}$. Therefore, according to Lemma~\ref{LDCHEchInt}, we have that $\mathrm{ip}^{(2,X)@}_{s}(\mathrm{CH}^{(2)}_{s}(\mathfrak{P}^{(2)}))$ is a second-order path of length strictly greater than one containing a  second-order echelon on its first step.

Moreover, taking into account Definition~\ref{DDIp}, we have that 
\begin{multline*}
\mathrm{ip}^{(2,X)@}_{s}\left(\mathrm{CH}^{(2)}_{s}\left(
\mathfrak{P}^{(2)}\right)\right)
\\
=
\mathrm{ip}^{(2,X)@}_{s}\left(\mathrm{CH}^{(2)}_{s}\left(
\mathfrak{P}^{(2),1,\bb{\mathfrak{P}^{(2)}}-1}
\right)\right)
\circ^{1\mathbf{Pth}_{\boldsymbol{\mathcal{A}}^{(2)}}}_{s}
\mathrm{ip}^{(2,X)@}_{s}\left(\mathrm{CH}^{(2)}_{s}\left(
\mathfrak{P}^{(2),0,0}
\right)\right).
\end{multline*}

By Proposition~\ref{PDIpDCH}, we have that 
\begin{enumerate}
\item $\mathrm{ip}^{(2,X)@}_{s}(\mathrm{CH}^{(2)}_{s}(\mathfrak{P}^{(2),0,0}))$ is a second-order path in $[\mathfrak{P}^{(2),0,0}]_{s}$. Therefore, according to Lemma~\ref{LDCHDEch}, we have that $\mathrm{ip}^{(2,X)@}_{s}(\mathrm{CH}^{(2)}_{s}(\mathfrak{P}^{(2),0,0}))$ is a  second-order echelon. 
\end{enumerate}

It follows by Definition~\ref{DDOrd} that the following two inequalities hold
\begin{align*}
&\left(
\mathrm{ip}^{(2,X)@}_{s}\left(
\mathrm{CH}^{(2)}_{s}\left(
\mathfrak{P}^{(2),0,0}
\right)\right),s\right)
\leq_{\mathbf{Pth}_{\boldsymbol{\mathcal{A}}^{(2)}}}
\left(
\mathrm{ip}^{(2,X)@}_{s}\left(
\mathrm{CH}^{(2)}_{s}\left(
\mathfrak{P}^{(2)}
\right)\right),s\right);
\\
&\left(
\mathrm{ip}^{(2,X)@}_{s}\left(
\mathrm{CH}^{(2)}_{s}\left(
\mathfrak{P}^{(2),1,\bb{\mathfrak{P}^{(2)}}-1}
\right)\right),s\right)
\leq_{\mathbf{Pth}_{\boldsymbol{\mathcal{A}}^{(2)}}}
\left(
\mathrm{ip}^{(2,X)@}_{s}\left(
\mathrm{CH}^{(2)}_{s}\left(
\mathfrak{P}^{(2)}
\right)\right),s\right).
\end{align*}

By assumption, we have that 
$
(\mathrm{CH}^{(2)}_{t}(\mathfrak{Q}^{(2)}),t)
\leq_{\mathbf{T}_{\Sigma^{\boldsymbol{\mathcal{A}}^{(2)}}}(X)}
(\mathrm{CH}^{(2)}_{s}(\mathfrak{P}^{(2)}),s)
$. Therefore, taking into account Assumption~\ref{ADIpDCHOrdI},  one of the following cases hold
\begin{enumerate}
\item[(1.1.1)] $s=t$ and $\mathrm{CH}^{(2)}_{t}(\mathfrak{Q}^{(2)})=\mathrm{CH}^{(2)}_{s}(\mathfrak{P}^{(2),0,0})$;
\item[(1.1.2)] $s=t$ and $\mathrm{CH}^{(2)}_{t}(\mathfrak{Q}^{(2)})=\mathrm{CH}^{(2)}_{s}(\mathfrak{P}^{(2),1,\bb{\mathfrak{P}^{(2)}}-1})$;
\item[(1.1.3)] $(\mathrm{CH}^{(2)}_{t}(\mathfrak{Q}^{(2)}),t)
\leq_{\mathbf{T}_{\Sigma^{\boldsymbol{\mathcal{A}}^{(2)}}}(X)}
(\mathrm{CH}^{(2)}_{s}(\mathfrak{P}^{(2),0,0}),s)$;
\item[(1.1.4)] $(\mathrm{CH}^{(2)}_{t}(\mathfrak{Q}^{(2)}),t)
\leq_{\mathbf{T}_{\Sigma^{\boldsymbol{\mathcal{A}}^{(2)}}}(X)}
(\mathrm{CH}^{(2)}_{s}(\mathfrak{P}^{(2),1,\bb{\mathfrak{P}^{(2)}}-1})$.
\end{enumerate}

If~(1.1.1), i.e., $t=s$ and $\mathrm{CH}^{(2)}_{t}(\mathfrak{Q}^{(2)})=\mathrm{CH}^{(2)}_{s}(\mathfrak{P}^{(2),0,0})$, we have already proven that
\[
\left(
\mathrm{ip}^{(2,X)@}_{s}\left(
\mathrm{CH}^{(2)}_{s}\left(
\mathfrak{P}^{(2),0,0}
\right)\right),s\right)
\leq_{\mathbf{Pth}_{\boldsymbol{\mathcal{A}}^{(2)}}}
\left(
\mathrm{ip}^{(2,X)@}_{s}\left(
\mathrm{CH}^{(2)}_{s}\left(
\mathfrak{P}^{(2)}
\right)\right),s\right).
\]

If~(1.1.2), i.e., $t=s$ and $\mathrm{CH}^{(2)}_{t}(\mathfrak{Q}^{(2)})=\mathrm{CH}^{(2)}_{s}(\mathfrak{P}^{(2),1,\bb{\mathfrak{P}^{(2)}}-1})$, we have already proven that
\[
\left(
\mathrm{ip}^{(2,X)@}_{s}\left(
\mathrm{CH}^{(2)}_{s}\left(
\mathfrak{P}^{(2),1,\bb{\mathfrak{P}^{(2)}}-1}
\right)\right),s\right)
\leq_{\mathbf{Pth}_{\boldsymbol{\mathcal{A}}^{(2)}}}
\left(
\mathrm{ip}^{(2,X)@}_{s}\left(
\mathrm{CH}^{(2)}_{s}\left(
\mathfrak{P}^{(2)}
\right)\right),s\right).
\]

If~(1.1.3), i.e., if $(\mathrm{CH}^{(2)}_{t}(\mathfrak{Q}^{(2)}),t)
\leq_{\mathbf{T}_{\Sigma^{\boldsymbol{\mathcal{A}}^{(2)}}}(X)}
(\mathrm{CH}^{(2)}_{s}(\mathfrak{P}^{(2),0,0}),s)$, then note that, according to Definition~\ref{DDOrd} $(\mathfrak{P}^{(2),0,0},s)\leq_{\mathbf{Pth}_{\boldsymbol{\mathcal{A}}^{(2)}}}(\mathfrak{P}^{(2)},s)$. Therefore, by induction we have that 
\[
\left(
\mathrm{ip}^{(2,X)@}_{t}\left(
\mathrm{CH}^{(2)}_{t}\left(
\mathfrak{Q}^{(2)}
\right)\right),t\right)
\leq_{\mathbf{Pth}_{\boldsymbol{\mathcal{A}}^{(2)}}}
\left(
\mathrm{ip}^{(2,X)@}_{s}\left(
\mathrm{CH}^{(2)}_{s}\left(
\mathfrak{P}^{(2),0,0}
\right)\right),s\right).
\]

Therefore, by transitivity of $
\leq_{\mathbf{Pth}_{\boldsymbol{\mathcal{A}}^{(2)}}}$, we conclude that 
\[
\left(
\mathrm{ip}^{(2,X)@}_{t}\left(
\mathrm{CH}^{(2)}_{t}\left(
\mathfrak{Q}^{(2)}
\right)\right),t\right)
\leq_{\mathbf{Pth}_{\boldsymbol{\mathcal{A}}^{(2)}}}
\left(
\mathrm{ip}^{(2,X)@}_{s}\left(
\mathrm{CH}^{(2)}_{s}\left(
\mathfrak{P}^{(2)}
\right)\right),s\right).
\]

If~(1.1.4), i.e., if $(\mathrm{CH}^{(2)}_{t}(\mathfrak{Q}^{(2)}),t)
\leq_{\mathbf{T}_{\Sigma^{\boldsymbol{\mathcal{A}}^{(2)}}}(X)}
(\mathrm{CH}^{(2)}_{s}(\mathfrak{P}^{(2),1,\bb{\mathfrak{P}^{(2)}}-1}),s)$, then note that, according to Definition~\ref{DDOrd} $(\mathfrak{P}^{(2),1,\bb{\mathfrak{P}^{(2)}}-1},s)\leq_{\mathbf{Pth}_{\boldsymbol{\mathcal{A}}^{(2)}}}(\mathfrak{P}^{(2)},s)$. Therefore, by induction we have that 
\[
\left(
\mathrm{ip}^{(2,X)@}_{t}\left(
\mathrm{CH}^{(2)}_{t}\left(
\mathfrak{Q}^{(2)}
\right)\right),t\right)
\leq_{\mathbf{Pth}_{\boldsymbol{\mathcal{A}}^{(2)}}}
\left(
\mathrm{ip}^{(2,X)@}_{s}\left(
\mathrm{CH}^{(2)}_{s}\left(
\mathfrak{P}^{(2),1,\bb{\mathfrak{P}^{(2)}}-1}
\right)\right),s\right).
\]

Therefore, by transitivity of $
\leq_{\mathbf{Pth}_{\boldsymbol{\mathcal{A}}^{(2)}}}$, we conclude that 
\[
\left(
\mathrm{ip}^{(2,X)@}_{t}\left(
\mathrm{CH}^{(2)}_{t}\left(
\mathfrak{Q}^{(2)}
\right)\right),t\right)
\leq_{\mathbf{Pth}_{\boldsymbol{\mathcal{A}}^{(2)}}}
\left(
\mathrm{ip}^{(2,X)@}_{s}\left(
\mathrm{CH}^{(2)}_{s}\left(
\mathfrak{P}^{(2)}
\right)\right),s\right).
\]

This completes the proof of Case~(1.1).

If~(1.2) $i>0$, i.e., if $\mathfrak{P}^{(2)}$ has its first  second-order echelon on a step different from the initial one, then according to Definition~\ref{DDCH}, we have that 
\[
\mathrm{CH}^{(2)}_{s}\left(\mathfrak{P}^{(2)}\right)
=
\mathrm{CH}^{(2)}_{s}\left(\mathfrak{P}^{(2),i,\bb{\mathfrak{P}^{(2)}}-1}\right)
\circ^{1\mathbf{T}_{\Sigma^{\boldsymbol{\mathcal{A}}^{(2)}}}(X)}_{s}
\mathrm{CH}^{(2)}_{s}\left(\mathfrak{P}^{(2),0,i-1}\right).
\]

Taking into account Proposition~\ref{PDIpDCH}, we have that $\mathrm{ip}^{(2,X)@}_{s}(\mathrm{CH}^{(2)}_{s}(\mathfrak{P}^{(2)}))$ is a second-order path in $[\mathfrak{P}^{(2)}]_{s}$. Therefore, according to Lemma~\ref{LDCHEchNInt}, we have that $\mathrm{ip}^{(2,X)@}_{s}(\mathrm{CH}^{(2)}_{s}(\mathfrak{P}^{(2)}))$ is a second-order path of length strictly greater than one containing its first  second-order echelon on a step different from the initial one.

Moreover, taking into account Definition~\ref{DDIp}, we have that 
\begin{multline*}
\mathrm{ip}^{(2,X)@}_{s}\left(\mathrm{CH}^{(2)}_{s}\left(
\mathfrak{P}^{(2)}\right)\right)
\\
=
\mathrm{ip}^{(2,X)@}_{s}\left(\mathrm{CH}^{(2)}_{s}\left(
\mathfrak{P}^{(2),i,\bb{\mathfrak{P}^{(2)}}-1}
\right)\right)
\circ^{1\mathbf{Pth}_{\boldsymbol{\mathcal{A}}^{(2)}}}_{s}
\mathrm{ip}^{(2,X)@}_{s}\left(\mathrm{CH}^{(2)}_{s}\left(
\mathfrak{P}^{(2),0,i-1}
\right)\right).
\end{multline*}

By Proposition~\ref{PDIpDCH}, we have that 
\begin{enumerate}
\item $\mathrm{ip}^{(2,X)@}_{s}(\mathrm{CH}^{(2)}_{s}(\mathfrak{P}^{(2),0,i-1}))$ is a second-order path in $[\mathfrak{P}^{(2),0,i-1}]_{s}$. Therefore, according to Lemma~\ref{LDCHNEch}, we have that $\mathrm{ip}^{(2,X)@}_{s}(\mathrm{CH}^{(2)}_{s}(\mathfrak{P}^{(2),0,i-1}))$ is an echelonless second-order path. 
\item $\mathrm{ip}^{(2,X)@}_{s}(\mathrm{CH}^{(2)}_{s}(\mathfrak{P}^{(2),i,\bb{\mathfrak{P}^{(2)}-1}}))$ is a second-order path in $[\mathfrak{P}^{(2),i,\bb{\mathfrak{P}^{(2)}}-1}]_{s}$. Therefore, according to either Lemmas~\ref{LDCHDEch} or~\ref{LDCHEchInt} , we have that $\mathrm{ip}^{(2,X)@}_{s}(\mathrm{CH}^{(2)}_{s}(\mathfrak{P}^{(2),i,\bb{\mathfrak{P}^{(2)}}-1}))$ is a second-order path containing a  second-order echelon on its first step. 
\end{enumerate}

It follows by Definition~\ref{DDOrd} that the following two inequalities hold
\[
\left(
\mathrm{ip}^{(2,X)@}_{s}\left(
\mathrm{CH}^{(2)}_{s}\left(
\mathfrak{P}^{(2),0,i-1}
\right)\right),s\right)
\leq_{\mathbf{Pth}_{\boldsymbol{\mathcal{A}}^{(2)}}}
\left(
\mathrm{ip}^{(2,X)@}_{s}\left(
\mathrm{CH}^{(2)}_{s}\left(
\mathfrak{P}^{(2)}
\right)\right),s\right);
\]
\[
\left(
\mathrm{ip}^{(2,X)@}_{s}\left(
\mathrm{CH}^{(2)}_{s}\left(
\mathfrak{P}^{(2),i,\bb{\mathfrak{P}^{(2)}}-1}
\right)\right),s\right)
\leq_{\mathbf{Pth}_{\boldsymbol{\mathcal{A}}^{(2)}}}
\left(
\mathrm{ip}^{(2,X)@}_{s}\left(
\mathrm{CH}^{(2)}_{s}\left(
\mathfrak{P}^{(2)}
\right)\right),s\right).
\]

By assumption, we have that 
$
(\mathrm{CH}^{(2)}_{t}(\mathfrak{Q}^{(2)}),t)
\leq_{\mathbf{T}_{\Sigma^{\boldsymbol{\mathcal{A}}^{(2)}}}(X)}
(\mathrm{CH}^{(2)}_{s}(\mathfrak{P}^{(2)}),s)
$. Therefore, taking into account Assumption~\ref{ADIpDCHOrdI},  one of the following cases hold
\begin{enumerate}
\item[(1.2.1)] $s=t$ and $\mathrm{CH}^{(2)}_{t}(\mathfrak{Q}^{(2)})=\mathrm{CH}^{(2)}_{s}(\mathfrak{P}^{(2),0,i-1})$;
\item[(1.2.2)] $s=t$ and $\mathrm{CH}^{(2)}_{t}(\mathfrak{Q}^{(2)})=\mathrm{CH}^{(2)}_{s}(\mathfrak{P}^{(2),i,\bb{\mathfrak{P}^{(2)}}-1})$;
\item[(1.2.3)] $(\mathrm{CH}^{(2)}_{t}(\mathfrak{Q}^{(2)}),t)
\leq_{\mathbf{T}_{\Sigma^{\boldsymbol{\mathcal{A}}^{(2)}}}(X)}
(\mathrm{CH}^{(2)}_{s}(\mathfrak{P}^{(2),0,i-1}),s)$;
\item[(1.2.4)] $(\mathrm{CH}^{(2)}_{t}(\mathfrak{Q}^{(2)}),t)
\leq_{\mathbf{T}_{\Sigma^{\boldsymbol{\mathcal{A}}^{(2)}}}(X)}
(\mathrm{CH}^{(2)}_{s}(\mathfrak{P}^{(2),i,\bb{\mathfrak{P}^{(2)}}-1})$.
\end{enumerate}

If~(1.2.1), i.e., $t=s$ and $\mathrm{CH}^{(2)}_{t}(\mathfrak{Q}^{(2)})=\mathrm{CH}^{(2)}_{s}(\mathfrak{P}^{(2),0,i-1})$, we have already proven that
\[
\left(
\mathrm{ip}^{(2,X)@}_{s}\left(
\mathrm{CH}^{(2)}_{s}\left(
\mathfrak{P}^{(2),0,i-1}
\right)\right),s\right)
\leq_{\mathbf{Pth}_{\boldsymbol{\mathcal{A}}^{(2)}}}
\left(
\mathrm{ip}^{(2,X)@}_{s}\left(
\mathrm{CH}^{(2)}_{s}\left(
\mathfrak{P}^{(2)}
\right)\right),s\right).
\]

If~(1.2.2), i.e., $t=s$ and $\mathrm{CH}^{(2)}_{t}(\mathfrak{Q}^{(2)})=\mathrm{CH}^{(2)}_{s}(\mathfrak{P}^{(2),i,\bb{\mathfrak{P}^{(2)}}-1})$, we have already proven that
\[
\left(
\mathrm{ip}^{(2,X)@}_{s}\left(
\mathrm{CH}^{(2)}_{s}\left(
\mathfrak{P}^{(2),i,\bb{\mathfrak{P}^{(2)}}-1}
\right)\right),s\right)
\leq_{\mathbf{Pth}_{\boldsymbol{\mathcal{A}}^{(2)}}}
\left(
\mathrm{ip}^{(2,X)@}_{s}\left(
\mathrm{CH}^{(2)}_{s}\left(
\mathfrak{P}^{(2)}
\right)\right),s\right).
\]

If~(1.2.3), i.e., if $(\mathrm{CH}^{(2)}_{t}(\mathfrak{Q}^{(2)}),t)
\leq_{\mathbf{T}_{\Sigma^{\boldsymbol{\mathcal{A}}^{(2)}}}(X)}
(\mathrm{CH}^{(2)}_{s}(\mathfrak{P}^{(2),0,i-1}),s)$, then note that, according to Definition~\ref{DDOrd} $(\mathfrak{P}^{(2),0,i-1},s)\leq_{\mathbf{Pth}_{\boldsymbol{\mathcal{A}}^{(2)}}}(\mathfrak{P}^{(2)},s)$. Therefore, by induction we have that 
\[
\left(
\mathrm{ip}^{(2,X)@}_{t}\left(
\mathrm{CH}^{(2)}_{t}\left(
\mathfrak{Q}^{(2)}
\right)\right),t\right)
\leq_{\mathbf{Pth}_{\boldsymbol{\mathcal{A}}^{(2)}}}
\left(
\mathrm{ip}^{(2,X)@}_{s}\left(
\mathrm{CH}^{(2)}_{s}\left(
\mathfrak{P}^{(2),0,i-1}
\right)\right),s\right).
\]

Therefore, by transitivity of $
\leq_{\mathbf{Pth}_{\boldsymbol{\mathcal{A}}^{(2)}}}$, we conclude that 
\[
\left(
\mathrm{ip}^{(2,X)@}_{t}\left(
\mathrm{CH}^{(2)}_{t}\left(
\mathfrak{Q}^{(2)}
\right)\right),t\right)
\leq_{\mathbf{Pth}_{\boldsymbol{\mathcal{A}}^{(2)}}}
\left(
\mathrm{ip}^{(2,X)@}_{s}\left(
\mathrm{CH}^{(2)}_{s}\left(
\mathfrak{P}^{(2)}
\right)\right),s\right).
\]

If~(1.2.4), i.e., if $(\mathrm{CH}^{(2)}_{t}(\mathfrak{Q}^{(2)}),t)
\leq_{\mathbf{T}_{\Sigma^{\boldsymbol{\mathcal{A}}^{(2)}}}(X)}
(\mathrm{CH}^{(2)}_{s}(\mathfrak{P}^{(2),i,\bb{\mathfrak{P}^{(2)}}-1}),s)$, then note that, according to Definition~\ref{DDOrd} $(\mathfrak{P}^{(2),i,\bb{\mathfrak{P}^{(2)}}-1},s)\leq_{\mathbf{Pth}_{\boldsymbol{\mathcal{A}}^{(2)}}}(\mathfrak{P}^{(2)},s)$. Therefore, by induction we have that 
\[
\left(
\mathrm{ip}^{(2,X)@}_{t}\left(
\mathrm{CH}^{(2)}_{t}\left(
\mathfrak{Q}^{(2)}
\right)\right),t\right)
\leq_{\mathbf{Pth}_{\boldsymbol{\mathcal{A}}^{(2)}}}
\left(
\mathrm{ip}^{(2,X)@}_{s}\left(
\mathrm{CH}^{(2)}_{s}\left(
\mathfrak{P}^{(2),i,\bb{\mathfrak{P}^{(2)}}-1}
\right)\right),s\right).
\]

Therefore, by transitivity of $
\leq_{\mathbf{Pth}_{\boldsymbol{\mathcal{A}}^{(2)}}}$, we conclude that 
\[
\left(
\mathrm{ip}^{(2,X)@}_{t}\left(
\mathrm{CH}^{(2)}_{t}\left(
\mathfrak{Q}^{(2)}
\right)\right),t\right)
\leq_{\mathbf{Pth}_{\boldsymbol{\mathcal{A}}^{(2)}}}
\left(
\mathrm{ip}^{(2,X)@}_{s}\left(
\mathrm{CH}^{(2)}_{s}\left(
\mathfrak{P}^{(2)}
\right)\right),s\right).
\]

This completes the proof of Case~(1.2).

This complete the proof of Case~(1).

If~(2), i.e., if $\mathfrak{P}^{(2)}$ is an echelonless second-order path, it could be the case that (2.1) $\mathfrak{P}^{(2)}$ is an echelonless second-order path that is not head constant, or (2.2) $\mathfrak{P}^{(2)}$ is a head-constant echelonless second-order path that is not coherent, or (2.3) $\mathfrak{P}^{(2)}$ is a coherent head-constant echelonless second-order path.

If~(2.1), i.e., if $\mathfrak{P}^{(2)}$ is an echelonless second-order path that is not head-constant, then we let $i\in\bb{\mathfrak{P}^{(2)}}$ be the greatest index for which $\mathfrak{P}^{(2),0,i}$ is a head-constant echelonless second-order path, then according to Definition~\ref{DDCH}, we have that 
\[
\mathrm{CH}^{(2)}_{s}\left(\mathfrak{P}^{(2)}\right)
=
\mathrm{CH}^{(2)}_{s}\left(\mathfrak{P}^{(2),i+1,\bb{\mathfrak{P}^{(2)}}-1}\right)
\circ^{1\mathbf{T}_{\Sigma^{\boldsymbol{\mathcal{A}}^{(2)}}}(X)}_{s}
\mathrm{CH}^{(2)}_{s}\left(\mathfrak{P}^{(2),0,i}\right).
\]

Taking into account Proposition~\ref{PDIpDCH}, we have that $\mathrm{ip}^{(2,X)@}_{s}(\mathrm{CH}^{(2)}_{s}(\mathfrak{P}^{(2)}))$ is a second-order path in $[\mathfrak{P}^{(2)}]_{s}$. Therefore, according to Lemma~\ref{LDCHNEchNHd}, we have that $\mathrm{ip}^{(2,X)@}_{s}(\mathrm{CH}^{(2)}_{s}(\mathfrak{P}^{(2)}))$ is an echelonless second-order path that is not head-constant.

Moreover, taking into account Definition~\ref{DDIp}, we have that 
\begin{multline*}
\mathrm{ip}^{(2,X)@}_{s}\left(\mathrm{CH}^{(2)}_{s}\left(
\mathfrak{P}^{(2)}\right)\right)
\\
=
\mathrm{ip}^{(2,X)@}_{s}\left(\mathrm{CH}^{(2)}_{s}\left(
\mathfrak{P}^{(2),i+1,\bb{\mathfrak{P}^{(2)}}-1}
\right)\right)
\circ^{1\mathbf{Pth}_{\boldsymbol{\mathcal{A}}^{(2)}}}_{s}
\mathrm{ip}^{(2,X)@}_{s}\left(\mathrm{CH}^{(2)}_{s}\left(
\mathfrak{P}^{(2),0,i}
\right)\right).
\end{multline*}

By Proposition~\ref{PDIpDCH}, we have that 
\begin{enumerate}
\item $\mathrm{ip}^{(2,X)@}_{s}(\mathrm{CH}^{(2)}_{s}(\mathfrak{P}^{(2),0,i}))$ is a second-order path in $[\mathfrak{P}^{(2),0,i}]_{s}$. Therefore, according to Lemma~\ref{LDCHNEchHd}, we have that $\mathrm{ip}^{(2,X)@}_{s}(\mathrm{CH}^{(2)}_{s}(\mathfrak{P}^{(2),0,i}))$ is a head-constant echelonless second-order path. Note that according to Lemma~\ref{LDCH}, we have that 
\[
\mathrm{tg}^{([1],2)}_{s}\left(
\mathrm{ip}^{(2,X)@}_{s}\left(
\mathrm{CH}^{(2)}_{s}\left(
\mathfrak{P}^{(2),0,i}
\right)\right)\right)
=\mathrm{tg}^{([1],2)}_{s}\left(
\mathfrak{P}^{(2),0,i}
\right).
\]
\item $\mathrm{ip}^{(2,X)@}_{s}(\mathrm{CH}^{(2)}_{s}(\mathfrak{P}^{(2),i+1,\bb{\mathfrak{P}^{(2)}-1}}))$ is a second-order path in $[\mathfrak{P}^{(2),i+1,\bb{\mathfrak{P}^{(2)}}-1}]_{s}$. Therefore, according to Lemma~\ref{LDCHNEch}, we have that $\mathrm{ip}^{(2,X)@}_{s}(\mathrm{CH}^{(2)}_{s}(\mathfrak{P}^{(2),i+1,\bb{\mathfrak{P}^{(2)}}-1}))$ is an echelonless second-order path. Note that according to Lemma~\ref{LDCH}, we have that 
\[
\mathrm{sc}^{([1],2)}_{s}\left(
\mathrm{ip}^{(2,X)@}_{s}\left(
\mathrm{CH}^{(2)}_{s}\left(
\mathfrak{P}^{(2),i+1,\bb{\mathfrak{P}^{(2)}-1}}
\right)\right)\right)
=\mathrm{sc}^{([1],2)}_{s}\left(
\mathfrak{P}^{(2),i+1,\bb{\mathfrak{P}^{(2)}-1}}
\right).
\]
\end{enumerate}

Therefore, $i\in\bb{\mathfrak{P}^{(2)}}$ is the greatest index for which the second-order path $(\mathrm{ip}^{(2,X)@}_{s}(\mathrm{CH}^{(2)}_{s}(\mathfrak{P}^{(2)})))^{0,i}$ is head-constant.  It follows by Definition~\ref{DDOrd} that the following two inequalities hold
\[\left(
\mathrm{ip}^{(2,X)@}_{s}\left(
\mathrm{CH}^{(2)}_{s}\left(
\mathfrak{P}^{(2),0,i}
\right)\right),s\right)
\leq_{\mathbf{Pth}_{\boldsymbol{\mathcal{A}}^{(2)}}}
\left(
\mathrm{ip}^{(2,X)@}_{s}\left(
\mathrm{CH}^{(2)}_{s}\left(
\mathfrak{P}^{(2)}
\right)\right),s\right);
\]
\[
\left(
\mathrm{ip}^{(2,X)@}_{s}\left(
\mathrm{CH}^{(2)}_{s}\left(
\mathfrak{P}^{(2),i+1,\bb{\mathfrak{P}^{(2)}}-1}
\right)\right),s\right)
\leq_{\mathbf{Pth}_{\boldsymbol{\mathcal{A}}^{(2)}}}
\left(
\mathrm{ip}^{(2,X)@}_{s}\left(
\mathrm{CH}^{(2)}_{s}\left(
\mathfrak{P}^{(2)}
\right)\right),s\right).
\]

By assumption, we have that 
$
(\mathrm{CH}^{(2)}_{t}(\mathfrak{Q}^{(2)}),t)
\leq_{\mathbf{T}_{\Sigma^{\boldsymbol{\mathcal{A}}^{(2)}}}(X)}
(\mathrm{CH}^{(2)}_{s}(\mathfrak{P}^{(2)}),s)
$. Therefore, taking into account Assumption~\ref{ADIpDCHOrdI},  one of the following cases hold
\begin{enumerate}
\item[(2.1.1)] $s=t$ and $\mathrm{CH}^{(2)}_{t}(\mathfrak{Q}^{(2)})=\mathrm{CH}^{(2)}_{s}(\mathfrak{P}^{(2),0,i})$;
\item[(2.1.2)] $s=t$ and $\mathrm{CH}^{(2)}_{t}(\mathfrak{Q}^{(2)})=\mathrm{CH}^{(2)}_{s}(\mathfrak{P}^{(2),i+1,\bb{\mathfrak{P}^{(2)}}-1})$;
\item[(2.1.3)] $(\mathrm{CH}^{(2)}_{t}(\mathfrak{Q}^{(2)}),t)
\leq_{\mathbf{T}_{\Sigma^{\boldsymbol{\mathcal{A}}^{(2)}}}(X)}
(\mathrm{CH}^{(2)}_{s}(\mathfrak{P}^{(2),0,i}),s)$;
\item[(2.1.4)] $(\mathrm{CH}^{(2)}_{t}(\mathfrak{Q}^{(2)}),t)
\leq_{\mathbf{T}_{\Sigma^{\boldsymbol{\mathcal{A}}^{(2)}}}(X)}
(\mathrm{CH}^{(2)}_{s}(\mathfrak{P}^{(2),i+,\bb{\mathfrak{P}^{(2)}}-1})$.
\end{enumerate}

If~(2.1.1), i.e., $t=s$ and $\mathrm{CH}^{(2)}_{t}(\mathfrak{Q}^{(2)})=\mathrm{CH}^{(2)}_{s}(\mathfrak{P}^{(2),0,i})$, we have already proven that
\[
\left(
\mathrm{ip}^{(2,X)@}_{s}\left(
\mathrm{CH}^{(2)}_{s}\left(
\mathfrak{P}^{(2),0,i}
\right)\right),s\right)
\leq_{\mathbf{Pth}_{\boldsymbol{\mathcal{A}}^{(2)}}}
\left(
\mathrm{ip}^{(2,X)@}_{s}\left(
\mathrm{CH}^{(2)}_{s}\left(
\mathfrak{P}^{(2)}
\right)\right),s\right).
\]

If~(2.1.2), i.e., $t=s$ and $\mathrm{CH}^{(2)}_{t}(\mathfrak{Q}^{(2)})=\mathrm{CH}^{(2)}_{s}(\mathfrak{P}^{(2),i+1,\bb{\mathfrak{P}^{(2)}}-1})$, we have already proven that
\[
\left(
\mathrm{ip}^{(2,X)@}_{s}\left(
\mathrm{CH}^{(2)}_{s}\left(
\mathfrak{P}^{(2),i+1,\bb{\mathfrak{P}^{(2)}}-1}
\right)\right),s\right)
\leq_{\mathbf{Pth}_{\boldsymbol{\mathcal{A}}^{(2)}}}
\left(
\mathrm{ip}^{(2,X)@}_{s}\left(
\mathrm{CH}^{(2)}_{s}\left(
\mathfrak{P}^{(2)}
\right)\right),s\right).
\]

If~(2.1.3), i.e., if $(\mathrm{CH}^{(2)}_{t}(\mathfrak{Q}^{(2)}),t)
\leq_{\mathbf{T}_{\Sigma^{\boldsymbol{\mathcal{A}}^{(2)}}}(X)}
(\mathrm{CH}^{(2)}_{s}(\mathfrak{P}^{(2),0,i}),s)$, then note that, according to Definition~\ref{DDOrd} $(\mathfrak{P}^{(2),0,i},s)\leq_{\mathbf{Pth}_{\boldsymbol{\mathcal{A}}^{(2)}}}(\mathfrak{P}^{(2)},s)$. Therefore, by induction we have that 
\[
\left(
\mathrm{ip}^{(2,X)@}_{t}\left(
\mathrm{CH}^{(2)}_{t}\left(
\mathfrak{Q}^{(2)}
\right)\right),t\right)
\leq_{\mathbf{Pth}_{\boldsymbol{\mathcal{A}}^{(2)}}}
\left(
\mathrm{ip}^{(2,X)@}_{s}\left(
\mathrm{CH}^{(2)}_{s}\left(
\mathfrak{P}^{(2),0,i}
\right)\right),s\right).
\]

Therefore, by transitivity of $
\leq_{\mathbf{Pth}_{\boldsymbol{\mathcal{A}}^{(2)}}}$, we conclude that 
\[
\left(
\mathrm{ip}^{(2,X)@}_{t}\left(
\mathrm{CH}^{(2)}_{t}\left(
\mathfrak{Q}^{(2)}
\right)\right),t\right)
\leq_{\mathbf{Pth}_{\boldsymbol{\mathcal{A}}^{(2)}}}
\left(
\mathrm{ip}^{(2,X)@}_{s}\left(
\mathrm{CH}^{(2)}_{s}\left(
\mathfrak{P}^{(2)}
\right)\right),s\right).
\]

If~(2.1.4), i.e., if $(\mathrm{CH}^{(2)}_{t}(\mathfrak{Q}^{(2)}),t)
\leq_{\mathbf{T}_{\Sigma^{\boldsymbol{\mathcal{A}}^{(2)}}}(X)}
(\mathrm{CH}^{(2)}_{s}(\mathfrak{P}^{(2),i+1,\bb{\mathfrak{P}^{(2)}}-1}),s)$, then note that, according to Definition~\ref{DDOrd} $(\mathfrak{P}^{(2),i+1,\bb{\mathfrak{P}^{(2)}}-1},s)\leq_{\mathbf{Pth}_{\boldsymbol{\mathcal{A}}^{(2)}}}(\mathfrak{P}^{(2)},s)$. Therefore, by induction we have that 
\[
\left(
\mathrm{ip}^{(2,X)@}_{t}\left(
\mathrm{CH}^{(2)}_{t}\left(
\mathfrak{Q}^{(2)}
\right)\right),t\right)
\leq_{\mathbf{Pth}_{\boldsymbol{\mathcal{A}}^{(2)}}}
\left(
\mathrm{ip}^{(2,X)@}_{s}\left(
\mathrm{CH}^{(2)}_{s}\left(
\mathfrak{P}^{(2),i+1,\bb{\mathfrak{P}^{(2)}}-1}
\right)\right),s\right).
\]

Therefore, by transitivity of $
\leq_{\mathbf{Pth}_{\boldsymbol{\mathcal{A}}^{(2)}}}$, we conclude that 
\[
\left(
\mathrm{ip}^{(2,X)@}_{t}\left(
\mathrm{CH}^{(2)}_{t}\left(
\mathfrak{Q}^{(2)}
\right)\right),t\right)
\leq_{\mathbf{Pth}_{\boldsymbol{\mathcal{A}}^{(2)}}}
\left(
\mathrm{ip}^{(2,X)@}_{s}\left(
\mathrm{CH}^{(2)}_{s}\left(
\mathfrak{P}^{(2)}
\right)\right),s\right).
\]

This completes the proof of Case~(2.1).

If~(2.2), i.e., if $\mathfrak{P}^{(2)}$ is a head-constant echelonless second-order path that is not coherent, then we let $i\in\bb{\mathfrak{P}^{(2)}}$ be the greatest index for which $\mathfrak{P}^{(2),0,i}$ is a coherent head-constant echelonless second-order path, then according to Definition~\ref{DDCH}, we have that 
\[
\mathrm{CH}^{(2)}_{s}\left(\mathfrak{P}^{(2)}\right)
=
\mathrm{CH}^{(2)}_{s}\left(\mathfrak{P}^{(2),i+1,\bb{\mathfrak{P}^{(2)}}-1}\right)
\circ^{1\mathbf{T}_{\Sigma^{\boldsymbol{\mathcal{A}}^{(2)}}}(X)}_{s}
\mathrm{CH}^{(2)}_{s}\left(\mathfrak{P}^{(2),0,i}\right).
\]

Taking into account Proposition~\ref{PDIpDCH}, we have that $\mathrm{ip}^{(2,X)@}_{s}(\mathrm{CH}^{(2)}_{s}(\mathfrak{P}^{(2)}))$ is a second-order path in $[\mathfrak{P}^{(2)}]_{s}$. Therefore, according to Lemma~\ref{LDCHNEchHdNC}, we have that $\mathrm{ip}^{(2,X)@}_{s}(\mathrm{CH}^{(2)}_{s}(\mathfrak{P}^{(2)}))$ is a head-constant echelonless second-order path that is not coherent.

Moreover, taking into account Definition~\ref{DDIp}, we have that 
\begin{multline*}
\mathrm{ip}^{(2,X)@}_{s}\left(\mathrm{CH}^{(2)}_{s}\left(
\mathfrak{P}^{(2)}\right)\right)
\\
=
\mathrm{ip}^{(2,X)@}_{s}\left(\mathrm{CH}^{(2)}_{s}\left(
\mathfrak{P}^{(2),i+1,\bb{\mathfrak{P}^{(2)}}-1}
\right)\right)
\circ^{1\mathbf{Pth}_{\boldsymbol{\mathcal{A}}^{(2)}}}_{s}
\mathrm{ip}^{(2,X)@}_{s}\left(\mathrm{CH}^{(2)}_{s}\left(
\mathfrak{P}^{(2),0,i}
\right)\right).
\end{multline*}

By Proposition~\ref{PDIpDCH}, we have that 
\begin{enumerate}
\item $\mathrm{ip}^{(2,X)@}_{s}(\mathrm{CH}^{(2)}_{s}(\mathfrak{P}^{(2),0,i}))$ is a second-order path in $[\mathfrak{P}^{(2),0,i}]_{s}$. Therefore, according to Lemma~\ref{LDCHNEchHdC}, we have that $\mathrm{ip}^{(2,X)@}_{s}(\mathrm{CH}^{(2)}_{s}(\mathfrak{P}^{(2),0,i}))$ is a coherent head-constant echelonless second-order path. 
\item $\mathrm{ip}^{(2,X)@}_{s}(\mathrm{CH}^{(2)}_{s}(\mathfrak{P}^{(2),i+1,\bb{\mathfrak{P}^{(2)}-1}}))$ is a second-order path in $[\mathfrak{P}^{(2),i+1,\bb{\mathfrak{P}^{(2)}}-1}]_{s}$. Therefore, according to Lemma~\ref{LDCHNEchHd}, we have that $\mathrm{ip}^{(2,X)@}_{s}(\mathrm{CH}^{(2)}_{s}(\mathfrak{P}^{(2),i+1,\bb{\mathfrak{P}^{(2)}}-1}))$ is a head-constant echelonless second-order path. 
\end{enumerate}

Following Lemma~\ref{LTech}, $i\in\bb{\mathfrak{P}^{(2)}}$ is the greatest index for which the second-order path $(\mathrm{ip}^{(2,X)@}_{s}(\mathrm{CH}^{(2)}_{s}(\mathfrak{P}^{(2)})))^{0,i}$ is coherent and head-constant.  It follows by Definition~\ref{DDOrd} that the following two inequalities hold
\[\left(
\mathrm{ip}^{(2,X)@}_{s}\left(
\mathrm{CH}^{(2)}_{s}\left(
\mathfrak{P}^{(2),0,i}
\right)\right),s\right)
\leq_{\mathbf{Pth}_{\boldsymbol{\mathcal{A}}^{(2)}}}
\left(
\mathrm{ip}^{(2,X)@}_{s}\left(
\mathrm{CH}^{(2)}_{s}\left(
\mathfrak{P}^{(2)}
\right)\right),s\right);
\]
\[
\left(
\mathrm{ip}^{(2,X)@}_{s}\left(
\mathrm{CH}^{(2)}_{s}\left(
\mathfrak{P}^{(2),i+1,\bb{\mathfrak{P}^{(2)}}-1}
\right)\right),s\right)
\leq_{\mathbf{Pth}_{\boldsymbol{\mathcal{A}}^{(2)}}}
\left(
\mathrm{ip}^{(2,X)@}_{s}\left(
\mathrm{CH}^{(2)}_{s}\left(
\mathfrak{P}^{(2)}
\right)\right),s\right).
\]

By assumption, we have that 
$
(\mathrm{CH}^{(2)}_{t}(\mathfrak{Q}^{(2)}),t)
\leq_{\mathbf{T}_{\Sigma^{\boldsymbol{\mathcal{A}}^{(2)}}}(X)}
(\mathrm{CH}^{(2)}_{s}(\mathfrak{P}^{(2)}),s)
$. Therefore, taking into account Assumption~\ref{ADIpDCHOrdI},  one of the following cases hold
\begin{enumerate}
\item[(2.2.1)] $s=t$ and $\mathrm{CH}^{(2)}_{t}(\mathfrak{Q}^{(2)})=\mathrm{CH}^{(2)}_{s}(\mathfrak{P}^{(2),0,i})$;
\item[(2.2.2)] $s=t$ and $\mathrm{CH}^{(2)}_{t}(\mathfrak{Q}^{(2)})=\mathrm{CH}^{(2)}_{s}(\mathfrak{P}^{(2),i+1,\bb{\mathfrak{P}^{(2)}}-1})$;
\item[(2.2.3)] $(\mathrm{CH}^{(2)}_{t}(\mathfrak{Q}^{(2)}),t)
\leq_{\mathbf{T}_{\Sigma^{\boldsymbol{\mathcal{A}}^{(2)}}}(X)}
(\mathrm{CH}^{(2)}_{s}(\mathfrak{P}^{(2),0,i}),s)$;
\item[(2.2.4)] $(\mathrm{CH}^{(2)}_{t}(\mathfrak{Q}^{(2)}),t)
\leq_{\mathbf{T}_{\Sigma^{\boldsymbol{\mathcal{A}}^{(2)}}}(X)}
(\mathrm{CH}^{(2)}_{s}(\mathfrak{P}^{(2),i+,\bb{\mathfrak{P}^{(2)}}-1})$.
\end{enumerate}

If~(2.2.1), i.e., $t=s$ and $\mathrm{CH}^{(2)}_{t}(\mathfrak{Q}^{(2)})=\mathrm{CH}^{(2)}_{s}(\mathfrak{P}^{(2),0,i})$, we have already proven that
\[
\left(
\mathrm{ip}^{(2,X)@}_{s}\left(
\mathrm{CH}^{(2)}_{s}\left(
\mathfrak{P}^{(2),0,i}
\right)\right),s\right)
\leq_{\mathbf{Pth}_{\boldsymbol{\mathcal{A}}^{(2)}}}
\left(
\mathrm{ip}^{(2,X)@}_{s}\left(
\mathrm{CH}^{(2)}_{s}\left(
\mathfrak{P}^{(2)}
\right)\right),s\right).
\]

If~(2.2.2), i.e., $t=s$ and $\mathrm{CH}^{(2)}_{t}(\mathfrak{Q}^{(2)})=\mathrm{CH}^{(2)}_{s}(\mathfrak{P}^{(2),i+1,\bb{\mathfrak{P}^{(2)}}-1})$, we have already proven that
\[
\left(
\mathrm{ip}^{(2,X)@}_{s}\left(
\mathrm{CH}^{(2)}_{s}\left(
\mathfrak{P}^{(2),i+1,\bb{\mathfrak{P}^{(2)}}-1}
\right)\right),s\right)
\leq_{\mathbf{Pth}_{\boldsymbol{\mathcal{A}}^{(2)}}}
\left(
\mathrm{ip}^{(2,X)@}_{s}\left(
\mathrm{CH}^{(2)}_{s}\left(
\mathfrak{P}^{(2)}
\right)\right),s\right).
\]

If~(2.2.3), i.e., if $(\mathrm{CH}^{(2)}_{t}(\mathfrak{Q}^{(2)}),t)
\leq_{\mathbf{T}_{\Sigma^{\boldsymbol{\mathcal{A}}^{(2)}}}(X)}
(\mathrm{CH}^{(2)}_{s}(\mathfrak{P}^{(2),0,i}),s)$, then note that, according to Definition~\ref{DDOrd} $(\mathfrak{P}^{(2),0,i},s)\leq_{\mathbf{Pth}_{\boldsymbol{\mathcal{A}}^{(2)}}}(\mathfrak{P}^{(2)},s)$. Therefore, by induction we have that 
\[
\left(
\mathrm{ip}^{(2,X)@}_{t}\left(
\mathrm{CH}^{(2)}_{t}\left(
\mathfrak{Q}^{(2)}
\right)\right),t\right)
\leq_{\mathbf{Pth}_{\boldsymbol{\mathcal{A}}^{(2)}}}
\left(
\mathrm{ip}^{(2,X)@}_{s}\left(
\mathrm{CH}^{(2)}_{s}\left(
\mathfrak{P}^{(2),0,i}
\right)\right),s\right).
\]

Therefore, by transitivity of $
\leq_{\mathbf{Pth}_{\boldsymbol{\mathcal{A}}^{(2)}}}$, we conclude that 
\[
\left(
\mathrm{ip}^{(2,X)@}_{t}\left(
\mathrm{CH}^{(2)}_{t}\left(
\mathfrak{Q}^{(2)}
\right)\right),t\right)
\leq_{\mathbf{Pth}_{\boldsymbol{\mathcal{A}}^{(2)}}}
\left(
\mathrm{ip}^{(2,X)@}_{s}\left(
\mathrm{CH}^{(2)}_{s}\left(
\mathfrak{P}^{(2)}
\right)\right),s\right).
\]

If~(2.2.4), i.e., if $(\mathrm{CH}^{(2)}_{t}(\mathfrak{Q}^{(2)}),t)
\leq_{\mathbf{T}_{\Sigma^{\boldsymbol{\mathcal{A}}^{(2)}}}(X)}
(\mathrm{CH}^{(2)}_{s}(\mathfrak{P}^{(2),i+1,\bb{\mathfrak{P}^{(2)}}-1}),s)$, then note that, according to Definition~\ref{DDOrd} $(\mathfrak{P}^{(2),i+1,\bb{\mathfrak{P}^{(2)}}-1},s)\leq_{\mathbf{Pth}_{\boldsymbol{\mathcal{A}}^{(2)}}}(\mathfrak{P}^{(2)},s)$. Therefore, by induction we have that 
\[
\left(
\mathrm{ip}^{(2,X)@}_{t}\left(
\mathrm{CH}^{(2)}_{t}\left(
\mathfrak{Q}^{(2)}
\right)\right),t\right)
\leq_{\mathbf{Pth}_{\boldsymbol{\mathcal{A}}^{(2)}}}
\left(
\mathrm{ip}^{(2,X)@}_{s}\left(
\mathrm{CH}^{(2)}_{s}\left(
\mathfrak{P}^{(2),i+1,\bb{\mathfrak{P}^{(2)}}-1}
\right)\right),s\right).
\]

Therefore, by transitivity of $
\leq_{\mathbf{Pth}_{\boldsymbol{\mathcal{A}}^{(2)}}}$, we conclude that 
\[
\left(
\mathrm{ip}^{(2,X)@}_{t}\left(
\mathrm{CH}^{(2)}_{t}\left(
\mathfrak{Q}^{(2)}
\right)\right),t\right)
\leq_{\mathbf{Pth}_{\boldsymbol{\mathcal{A}}^{(2)}}}
\left(
\mathrm{ip}^{(2,X)@}_{s}\left(
\mathrm{CH}^{(2)}_{s}\left(
\mathfrak{P}^{(2)}
\right)\right),s\right).
\]

This completes the proof of Case~(2.2).

If~(2.3), i.e., if $\mathfrak{P}^{(2)}$ is a coherent head-constant echelonless second-order path, then in virtue of Definition~\ref{DDHeadCt} there exists a unique word $\mathbf{s}\in S^{\star}-\{\lambda\}$ and a unique operation symbol $\tau\in\Sigma^{\boldsymbol{\mathcal{A}}}_{\mathbf{s},s}$ associated to $\mathfrak{P}^{(2)}$. Let $(\mathfrak{P}^{(2)}_{j})_{j\in\bb{\mathbf{s}}}$ be the family of second-order paths in $\mathrm{Pth}_{\boldsymbol{\mathcal{A}}^{(2)},\mathbf{s}}$ which, in virtue of Lemma~\ref{LDPthExtract} we can extract from $\mathfrak{P}^{(2)}$. Then according to Definition~\ref{DDCH}, we have that 
\[
\mathrm{CH}^{(2)}_{s}\left(\mathfrak{P}^{(2)}\right)
=
\tau^{\mathbf{T}_{\Sigma^{\boldsymbol{\mathcal{A}}^{(2)}}}(X)}\left(\left(
\mathrm{CH}^{(2)}_{s_{j}}\left(\mathfrak{P}^{(2)}_{j}\right)
\right)_{j\in\bb{\mathbf{s}}}\right)
.
\]

Taking into account Proposition~\ref{PDIpDCH}, we have that $\mathrm{ip}^{(2,X)@}_{s}(\mathrm{CH}^{(2)}_{s}(\mathfrak{P}^{(2)}))$ is a second-order path in $[\mathfrak{P}^{(2)}]_{s}$. Therefore, according to Lemma~\ref{LDCHNEchHdC}, we have that $\mathrm{ip}^{(2,X)@}_{s}(\mathrm{CH}^{(2)}_{s}(\mathfrak{P}^{(2)}))$ is a coherent head-constant echelonless second-order path associated to the same operation symbol $\tau\in\Sigma^{\boldsymbol{\mathcal{A}}}_{\mathbf{s},s}$.

Moreover, taking into account Definition~\ref{DDIp}, we have that 
\[
\mathrm{ip}^{(2,X)@}_{s}\left(\mathrm{CH}^{(2)}_{s}\left(
\mathfrak{P}^{(2)}\right)\right)
=
\tau^{\mathbf{Pth}_{\boldsymbol{\mathcal{A}}^{(2)}}}
\left(\left(
\mathrm{ip}^{(2,X)@}_{s_{j}}\left(\mathrm{CH}^{(2)}_{s_{j}}\left(
\mathfrak{P}^{(2)}_{j}
\right)\right)
\right)_{j\in\bb{\mathbf{s}}}\right)
.
\]

By Proposition~\ref{PDIpDCH}, we have that 
\begin{enumerate}
\item for every $j\in\bb{\mathbf{s}}$, $\mathrm{ip}^{(2,X)@}_{s_{j}}(\mathrm{CH}^{(2)}_{s_{j}}(\mathfrak{P}^{(2)}_{j}))$ is a second-order path in $[\mathfrak{P}^{(2)}_{j}]_{s_{j}}$. Moreover, the family of second-order paths that we can extract from $\mathrm{ip}^{(2,X)@}_{s}(\mathrm{CH}^{(2)}_{s}(\mathfrak{P}^{(2)}))$ in virtue of Lemma~\ref{LDPthExtract} is, precisely,
\[
\left(
\mathrm{ip}^{(2,X)@}_{s_{j}}\left(
\mathrm{CH}^{(2)}_{s_{j}}\left(
\mathfrak{P}^{(2)}_{j}
\right)\right)\right)_{j\in\bb{\mathbf{s}}}.
\]
\end{enumerate}

It follows by Definition~\ref{DDOrd} that, for every $j\in\bb{\mathbf{s}}$, the following inequality holds
\begin{align*}
&\left(
\mathrm{ip}^{(2,X)@}_{s_{j}}\left(
\mathrm{CH}^{(2)}_{s_{j}}\left(
\mathfrak{P}^{(2)}_{j}
\right)\right),s_{j}\right)
\leq_{\mathbf{Pth}_{\boldsymbol{\mathcal{A}}^{(2)}}}
\left(
\mathrm{ip}^{(2,X)@}_{s}\left(
\mathrm{CH}^{(2)}_{s}\left(
\mathfrak{P}^{(2)}
\right)\right),s\right).
\end{align*}

By assumption, we have that 
$
(\mathrm{CH}^{(2)}_{t}(\mathfrak{Q}^{(2)}),t)
\leq_{\mathbf{T}_{\Sigma^{\boldsymbol{\mathcal{A}}^{(2)}}}(X)}
(\mathrm{CH}^{(2)}_{s}(\mathfrak{P}^{(2)}),s)
$. Therefore, taking into account Assumption~\ref{ADIpDCHOrdI},  one of the following cases hold
\begin{enumerate}
\item[(2.3.1)] For some $j\in\bb{\mathbf{s}}$, $s_{j}=t$ and $\mathrm{CH}^{(2)}_{t}(\mathfrak{Q}^{(2)})=\mathrm{CH}^{(2)}_{s_{j}}(\mathfrak{P}^{(2)}_{j})$;
\item[(2.3.2)] For some $j\in\bb{\mathbf{s}}$, $(\mathrm{CH}^{(2)}_{t}(\mathfrak{Q}^{(2)}),t)
\leq_{\mathbf{T}_{\Sigma^{\boldsymbol{\mathcal{A}}^{(2)}}}(X)}
(\mathrm{CH}^{(2)}_{s}(\mathfrak{P}^{(2)}_{j}),s_{j})$.
\end{enumerate}

If~(2.3.1), i.e., $t=s_{j}$, for some $j\in\bb{\mathbf{s}}$, and $\mathrm{CH}^{(2)}_{t}(\mathfrak{Q}^{(2)})=\mathrm{CH}^{(2)}_{s_{j}}(\mathfrak{P}^{(2)}_{j})$, we have already proven that
\[
\left(
\mathrm{ip}^{(2,X)@}_{s_{j}}\left(
\mathrm{CH}^{(2)}_{s_{j}}\left(
\mathfrak{P}^{(2)}_{j}
\right)\right),s_{j}\right)
\leq_{\mathbf{Pth}_{\boldsymbol{\mathcal{A}}^{(2)}}}
\left(
\mathrm{ip}^{(2,X)@}_{s}\left(
\mathrm{CH}^{(2)}_{s}\left(
\mathfrak{P}^{(2)}
\right)\right),s\right).
\]

If~(2.3.2), i.e., if $(\mathrm{CH}^{(2)}_{t}(\mathfrak{Q}^{(2)}),t)
\leq_{\mathbf{T}_{\Sigma^{\boldsymbol{\mathcal{A}}^{(2)}}}(X)}
(\mathrm{CH}^{(2)}_{s_{j}}(\mathfrak{P}^{(2)}_{j}),s_{j})$, then note that, according to Definition~\ref{DDOrd} $(\mathfrak{P}^{(2)}_{j},s_{j})\leq_{\mathbf{Pth}_{\boldsymbol{\mathcal{A}}^{(2)}}}(\mathfrak{P}^{(2)},s)$. Therefore, by induction we have that 
\[
\left(
\mathrm{ip}^{(2,X)@}_{t}\left(
\mathrm{CH}^{(2)}_{t}\left(
\mathfrak{Q}^{(2)}
\right)\right),t\right)
\leq_{\mathbf{Pth}_{\boldsymbol{\mathcal{A}}^{(2)}}}
\left(
\mathrm{ip}^{(2,X)@}_{s_{j}}\left(
\mathrm{CH}^{(2)}_{s_{j}}\left(
\mathfrak{P}^{(2)}_{j}
\right)\right),s_{j}\right).
\]

Therefore, by transitivity of $
\leq_{\mathbf{Pth}_{\boldsymbol{\mathcal{A}}^{(2)}}}$, we conclude that 
\[
\left(
\mathrm{ip}^{(2,X)@}_{t}\left(
\mathrm{CH}^{(2)}_{t}\left(
\mathfrak{Q}^{(2)}
\right)\right),t\right)
\leq_{\mathbf{Pth}_{\boldsymbol{\mathcal{A}}^{(2)}}}
\left(
\mathrm{ip}^{(2,X)@}_{s}\left(
\mathrm{CH}^{(2)}_{s}\left(
\mathfrak{P}^{(2)}
\right)\right),s\right).
\]

This completes the proof of Case~(2.3).

This completes the proof.
\end{proof}

\begin{restatable}{corollary}{CDIpDCHOrd}
\label{CDIpDCHOrd} The mapping $\coprod\mathrm{ip}^{(2,X)@}\circ \mathrm{CH}^{(2)}$ from $\coprod\mathrm{Pth}_{\boldsymbol{\mathcal{A}}^{(2)}}$ to $\coprod\mathrm{Pth}_{\boldsymbol{\mathcal{A}}^{(2)}}$ is order-preserving
\[
\textstyle
\coprod\mathrm{ip}^{(2,X)@}\circ \mathrm{CH}^{(2)}
\colon
\left(
\coprod\mathrm{Pth}_{\boldsymbol{\mathcal{A}}^{(2)}},
\leq_{\mathbf{Pth}_{\boldsymbol{\mathcal{A}}^{(2)}}}
\right)
\mor
\left(
\coprod\mathrm{Pth}_{\boldsymbol{\mathcal{A}}^{(2)}},
\leq_{\mathbf{Pth}_{\boldsymbol{\mathcal{A}}^{(2)}}}
\right),
\]
that is, given pairs $(\mathfrak{Q}^{(2)},t)$, $(\mathfrak{P}^{(2)},s)$ in $\coprod\mathrm{Pth}_{\boldsymbol{\mathcal{A}}^{(2)}}$, if $(\mathfrak{Q}^{(2)},t)\leq_{\mathbf{Pth}_{\boldsymbol{\mathcal{A}}^{(2)}}}(\mathfrak{P}^{(2)},s)$ then
\[
\left(
\mathrm{ip}^{(2,X)@}_{t}\left(
\mathrm{CH}^{(2)}_{t}\left(
\mathfrak{Q}^{(2)}
\right)\right),t\right)
\leq_{\mathbf{Pth}_{\boldsymbol{\mathcal{A}}^{(2)}}}
\left(
\mathrm{ip}^{(2,X)@}_{s}\left(
\mathrm{CH}^{(2)}_{s}\left(
\mathfrak{P}^{(2)}
\right)\right),s\right).
\]
\end{restatable}
\begin{proof}
It follows from Proposition~\ref{PDCHMono} and Proposition~\ref{PDIpDCHOrd}.
\end{proof}

\chapter{
\texorpdfstring
{The correspondence between second-order paths and terms}
{The correspondence between second-order paths and terms}
}\label{S2J}

In this chapter we introduce the relation $\Theta^{(2)}$ in $\mathrm{T}_{\Sigma^{\boldsymbol{\mathcal{A}}^{(2)}}}(X)$, with the aim of matching different terms that, through $\mathrm{ip}^{(2,X)@}$, are sent to second-order paths in the same class for the kernel of the second-order Curry-Howard mapping. We first show that $\Theta^{(1)}$, when lifted via $\eta^{(2,1)\sharp}$, is included in  $\Theta^{(2)}$. We also consider the smallest $\Sigma^{\boldsymbol{\mathcal{A}}^{(2)}}$-congruence containing $\Theta^{(2)}$, denoted by $\Theta^{[2]}$, and we check that this congruence also maintains the property of matching  different terms that, through $\mathrm{ip}^{(2,X)@}$, are sent to second-order paths in the same class for the kernel of the second-order Curry-Howard mapping. We also show that $\Theta^{[1]}$, when lifted via $\eta^{(2,1)\sharp}$, is included in  $\Theta^{[2]}$. We can then consider the quotient $\Sigma^{\boldsymbol{\mathcal{A}}^{(2)}}$-algebra $\mathrm{T}_{\Sigma^{\boldsymbol{\mathcal{A}}^{(2)}}}(X)/{\Theta^{[2]}}$. We also show that $\Theta^{[1]}$, when lifted via $\eta^{(2,1)\sharp}$, is included in  $\Theta^{[2]}$. It is shown that the composition $\mathrm{pr}^{\Theta^{[2]}}\circ \mathrm{CH}^{(2)}$, of the second-order Curry-Howard mapping and the projection to the $\Theta^{[2]}$-class, is itself a $\Sigma^{\boldsymbol{\mathcal{A}}^{(2)}}$-homomorphism. Moreover, if a term $P$ in $\mathrm{T}_{\Sigma^{\boldsymbol{\mathcal{A}}^{(2)}}}(X)_{s}$ is such that $\mathrm{ip}^{(2,X)@}_{s}(P)$ is a second-order path in $\mathrm{Pth}_{\boldsymbol{\mathcal{A}}^{(2)},s}$, then the terms $P$ and $\mathrm{CH}^{(2)}_{s}(\mathrm{ip}^{(2,X)@}_{s}(P))$ are $\Theta^{[2]}_{s}$-related. Indeed, for a pair of terms $(P,Q)$ in $\Theta^{[2]}_{s}$, it is shown that $\mathrm{ip}^{(2,X)@}_{s}(P)$ is a second-order path in $\mathrm{Pth}_{\boldsymbol{\mathcal{A}}^{(2)},s}$ if, and only if, $\mathrm{ip}^{(2,X)@}_{s}(Q)$ is a second-order path in $\mathrm{Pth}_{\boldsymbol{\mathcal{A}}^{(2)},s}$. Moreover, when one of such conditions holds, the two second-order paths $\mathrm{ip}^{(2,X)@}_{s}(P)$ and $\mathrm{ip}^{(2,X)@}_{s}(Q)$ have the same image under the second-order Curry-Howard mapping. This, in particular implies that, if a term $P$ is $\Theta^{[2]}_{s}$-related with $\mathrm{CH}^{(2)}_{s}(\mathfrak{P}^{(2)})$, then $\mathrm{ip}^{(2,X)@}_{s}(P)$  is itself a second-order path in the equivalence class $[\mathfrak{P}^{(2)}]_{s}$. Furthermore it is shown that, for two second-order paths $\mathfrak{P}^{(2)}$ and $\mathfrak{P}'^{(2)}$ in $\mathrm{Pth}_{\boldsymbol{\mathcal{A}}^{(2)},s}$, if $\mathrm{CH}^{(2)}_{s}(\mathfrak{P}^{(2)})$ and $\mathrm{CH}^{(2)}_{s}(\mathfrak{P}'^{(2)})$ are $\Theta^{[2]}_{s}$-related then the equivalence classes $[\mathfrak{P}^{(2)}]_{s}$ and $[\mathfrak{P}'^{(2)}]_{s}$ are equal.


We begin by defining a relation $\Theta^{(2)}$ on the underlying $S$-sorted set of $\mathbf{T}_{\Sigma^{\boldsymbol{\mathcal{A}}^{(2)}}}(X)$, the free $\Sigma^{\boldsymbol{\mathcal{A}}^{(2)}}$-algebra on $X$ that intertwines $0$-sources, $0$-targets, $0$-composition, $1$-sources, $1$-targets, and $1$-composition with the second-order Curry-Howard mapping. Such a relation, as well as other relations associated with it, will play a great role in the subsequent work.

\begin{restatable}{definition}{DDTheta}
\label{DDTheta} 
\index{Theta!second-order!$\Theta^{(2)}$}
We define $\Theta^{(2)}=(\Theta^{(2)}_{s})_{s\in S}$ to be the relation on $\mathrm{T}_{\Sigma^{\boldsymbol{\mathcal{A}}^{(2)}}}(X)$ consisting exactly of the following pairs of terms 
\begin{itemize}

\item[(i)] For every sort $s\in S$ and every second-order path $\mathfrak{P}^{(2)}$ in $\mathrm{Pth}_{\boldsymbol{\mathcal{A}}^{(2)},s}$, 
$$
\left(\mathrm{CH}^{(2)}_{s}\left(
\mathrm{sc}^{0\mathbf{Pth}_{\boldsymbol{\mathcal{A}}^{(2)}}}_{s}\left(
\mathfrak{P}^{(2)}
\right)\right),
\mathrm{sc}^{0\mathbf{T}_{\Sigma^{\boldsymbol{\mathcal{A}}^{(2)}}}(X)}_{s}\left(
\mathrm{CH}^{(2)}_{s}\left(
\mathfrak{P}^{(2)}
\right)\right)
\right)\in\Theta^{(2)}_{s};
$$ 
\item[(ii)] For every sort $s\in S$ and every second-order path $\mathfrak{P}^{(2)}$ in $\mathrm{Pth}_{\boldsymbol{\mathcal{A}}^{(2)},s}$, 
$$
\left(\mathrm{CH}^{(2)}_{s}\left(
\mathrm{tg}^{0\mathbf{Pth}_{\boldsymbol{\mathcal{A}}^{(2)}}}_{s}\left(
\mathfrak{P}^{(2)}
\right)\right),
\mathrm{tg}^{0\mathbf{T}_{\Sigma^{\boldsymbol{\mathcal{A}}^{(2)}}}(X)}_{s}\left(
\mathrm{CH}^{(2)}_{s}\left(
\mathfrak{P}^{(2)}
\right)\right)
\right)\in\Theta^{(2)}_{s};
$$ 
\item[(iii)] For every sort $s\in S$ and every pair of second-order paths $\mathfrak{Q}^{(2)},\mathfrak{P}^{(2)}$ in $\mathrm{Pth}_{\boldsymbol{\mathcal{A}}^{(2)},s}$, if
$\mathrm{sc}^{(0,2)}_{s}(\mathfrak{Q}^{(2)})=\mathrm{tg}^{(0,2)}_{s}(\mathfrak{P}^{(2)})$, then
$$
\left(\mathrm{CH}^{(2)}_{s}\left(
\mathfrak{Q}^{(2)}
\circ^{0\mathbf{Pth}_{\boldsymbol{\mathcal{A}}^{(2)}}}_{s}
\mathfrak{P}^{(2)}
\right),
\mathrm{CH}^{(2)}_{s}\left(
\mathfrak{Q}^{(2)}
\right)
\circ^{0\mathbf{T}_{\Sigma^{\boldsymbol{\mathcal{A}}^{(2)}}}(X)}_{s}
\mathrm{CH}^{(2)}_{s}\left(
\mathfrak{P}^{(2)}
\right)
\right)\in\Theta^{(2)}_{s}.
$$ 
\item[(iv)] For every sort $s\in S$ and every second-order path $\mathfrak{P}^{(2)}$ in $\mathrm{Pth}_{\boldsymbol{\mathcal{A}}^{(2)},s}$, 
$$
\left(\mathrm{CH}^{(2)}_{s}\left(
\mathrm{sc}^{1\mathbf{Pth}_{\boldsymbol{\mathcal{A}}^{(2)}}}_{s}\left(
\mathfrak{P}^{(2)}
\right)\right),
\mathrm{sc}^{1\mathbf{T}_{\Sigma^{\boldsymbol{\mathcal{A}}^{(2)}}}(X)}_{s}\left(
\mathrm{CH}^{(2)}_{s}\left(
\mathfrak{P}^{(2)}
\right)\right)
\right)\in\Theta^{(2)}_{s};
$$ 
\item[(v)] For every sort $s\in S$ and every second-order path $\mathfrak{P}^{(2)}$ in $\mathrm{Pth}_{\boldsymbol{\mathcal{A}}^{(2)},s}$, 
$$
\left(\mathrm{CH}^{(2)}_{s}\left(
\mathrm{tg}^{1\mathbf{Pth}_{\boldsymbol{\mathcal{A}}^{(2)}}}_{s}\left(
\mathfrak{P}^{(2)}
\right)\right),
\mathrm{tg}^{1\mathbf{T}_{\Sigma^{\boldsymbol{\mathcal{A}}^{(2)}}}(X)}_{s}\left(
\mathrm{CH}^{(2)}_{s}\left(
\mathfrak{P}^{(2)}
\right)\right)
\right)\in\Theta^{(2)}_{s};
$$ 
\item[(vi)] For every sort $s\in S$ and every pair of second-order paths $\mathfrak{Q}^{(2)},\mathfrak{P}^{(2)}$ in $\mathrm{Pth}_{\boldsymbol{\mathcal{A}}^{(2)},s}$, if
$\mathrm{sc}^{([1],2)}_{s}(\mathfrak{Q}^{(2)})=\mathrm{tg}^{([1],2)}_{s}(\mathfrak{P}^{(2)})$, then
$$
\left(\mathrm{CH}^{(2)}_{s}\left(
\mathfrak{Q}^{(2)}
\circ^{1\mathbf{Pth}_{\boldsymbol{\mathcal{A}}^{(2)}}}_{s}
\mathfrak{P}^{(2)}
\right),
\mathrm{CH}^{(2)}_{s}\left(
\mathfrak{Q}^{(2)}
\right)
\circ^{1\mathbf{T}_{\Sigma^{\boldsymbol{\mathcal{A}}^{(2)}}}(X)}_{s}
\mathrm{CH}^{(2)}_{s}\left(
\mathfrak{P}^{(2)}
\right)
\right)\in\Theta^{(2)}_{s}.
$$ 
\end{itemize}
This completes the definition of $\Theta^{(2)}$.
\end{restatable}

%
%
%

We first note the similarities between this relation and the relation $\Theta^{(1)}$ defined on $\mathrm{T}_{\Sigma^{\boldsymbol{\mathcal{A}}}}(X)$ and introduced in Definition~\ref{DTheta}. On the previous occasion we asked for compatibility with the categorial operations in $\Sigma^{\boldsymbol{\mathcal{A}}}$, while now we will ask for compatibility with the categorial operations in $\Sigma^{\boldsymbol{\mathcal{A}}^{(2)}}$.

Below we see how the relation $\Theta^{(1)}$ is related to the relation $\Theta^{(2)}$.

\begin{restatable}{proposition}{PDThetaDU}
\label{PDThetaDU} For the relation $\Theta^{(1)}\subseteq \mathrm{T}_{\Sigma^{\boldsymbol{\mathcal{A}}}}(X)^{2}$ and the relation $\Theta^{(2)}\subseteq \mathrm{T}_{\Sigma^{\boldsymbol{\mathcal{A}}^{(2)}}}(X)^{2}$, the following inclusion holds
\[
\eta^{(2,1)\sharp}\times \eta^{(2,1)\sharp}\left[\Theta^{(1)}\right]\subseteq \Theta^{(2)}.
\]
\end{restatable}
\begin{proof}
Let us recall from Definition~\ref{DTheta} that the relation $\Theta^{(1)}\subseteq \mathrm{T}_{\Sigma^{\boldsymbol{\mathcal{A}}}}(X)^{2}$ contains the following pairs
\item[(i)] For every sort $s\in S$ and every  path $\mathfrak{P}$ in $\mathrm{Pth}_{\boldsymbol{\mathcal{A}},s}$, 
$$
\left(\mathrm{CH}^{(1)}_{s}\left(
\mathrm{sc}^{0\mathbf{Pth}_{\boldsymbol{\mathcal{A}}}}_{s}\left(
\mathfrak{P}
\right)\right),
\mathrm{sc}^{0\mathbf{T}_{\Sigma^{\boldsymbol{\mathcal{A}}}}(X)}_{s}\left(
\mathrm{CH}^{(1)}_{s}\left(
\mathfrak{P}
\right)\right)
\right)\in\Theta^{(1)}_{s};
$$ 
\item[(ii)] For every sort $s\in S$ and every  path $\mathfrak{P}$ in $\mathrm{Pth}_{\boldsymbol{\mathcal{A}},s}$, 
$$
\left(\mathrm{CH}^{(1)}_{s}\left(
\mathrm{tg}^{0\mathbf{Pth}_{\boldsymbol{\mathcal{A}}}}_{s}\left(
\mathfrak{P}
\right)\right),
\mathrm{tg}^{0\mathbf{T}_{\Sigma^{\boldsymbol{\mathcal{A}}}}(X)}_{s}\left(
\mathrm{CH}^{(1)}_{s}\left(
\mathfrak{P}
\right)\right)
\right)\in\Theta^{(1)}_{s};
$$ 
\item[(iii)] For every sort $s\in S$ and every pair of paths $\mathfrak{Q},\mathfrak{P}$ in $\mathrm{Pth}_{\boldsymbol{\mathcal{A}},s}$, if
$\mathrm{sc}^{(0,1)}_{s}(\mathfrak{Q})=\mathrm{tg}^{(0,1)}_{s}(\mathfrak{P})$, then
$$
\left(\mathrm{CH}^{(1)}_{s}\left(
\mathfrak{Q}
\circ^{0\mathbf{Pth}_{\boldsymbol{\mathcal{A}}}}_{s}
\mathfrak{P}
\right),
\mathrm{CH}^{(1)}_{s}\left(
\mathfrak{Q}
\right)
\circ^{0\mathbf{T}_{\Sigma^{\boldsymbol{\mathcal{A}}}}(X)}_{s}
\mathrm{CH}^{(1)}_{s}\left(
\mathfrak{P}
\right)
\right)\in\Theta^{(1)}_{s}.
$$ 

We will check that every case satisfies the desired property.

\textsf{Case~(1)} Let $s$ be a sort in $S$ and $\mathfrak{P}$ a path in $\mathrm{Pth}_{\boldsymbol{\mathcal{A}},s}$, we want to prove that 
\[
\left(
\eta^{(2,1)\sharp}_{s}\left(
\mathrm{CH}^{(1)}_{s}\left(
\mathrm{sc}^{0\mathbf{Pth}_{\boldsymbol{\mathcal{A}}}}_{s}\left(
\mathfrak{P}
\right)\right)
\right)
,
\eta^{(2,1)\sharp}_{s}\left(
\mathrm{sc}^{0\mathbf{T}_{\Sigma^{\boldsymbol{\mathcal{A}}}}(X)}_{s}\left(
\mathrm{CH}^{(1)}_{s}\left(
\mathfrak{P}
\right)\right)
\right)
\right)
\in\Theta^{(2)}_{s}.
\]

For this, consider the second-order path $\mathfrak{P}^{(2)}=\mathrm{ip}^{(2,[1])\sharp}_{s}([\mathrm{CH}^{(1)}_{s}(\mathfrak{P})]_{s})$, according to Definition~\ref{DDTheta}, the following pair is in $\Theta^{(2)}_{s}$, 
\begin{equation}
\left(\mathrm{CH}^{(2)}_{s}\left(
\mathrm{sc}^{0\mathbf{Pth}_{\boldsymbol{\mathcal{A}}^{(2)}}}_{s}\left(
\mathfrak{P}^{(2)}
\right)\right),
\mathrm{sc}^{0\mathbf{T}_{\Sigma^{\boldsymbol{\mathcal{A}}^{(2)}}}(X)}_{s}\left(
\mathrm{CH}^{(2)}_{s}\left(
\mathfrak{P}^{(2)}
\right)\right)
\right)\in\Theta^{(2)}_{s}.
\tag{E1}\label{EDThetaDUI}
\end{equation}

Consider the left hand side of the pair in Equation~\ref{EDThetaDUI}. The following chain of equalities holds
\begin{flushleft}
$\mathrm{CH}^{(2)}_{s}\left(
\mathrm{sc}^{0\mathbf{Pth}_{\boldsymbol{\mathcal{A}}^{(2)}}}_{s}\left(
\mathfrak{P}^{(2)}
\right)\right)$
\allowdisplaybreaks
\begin{align*}
&=
\mathrm{CH}^{(2)}_{s}\left(
\mathrm{sc}^{0\mathbf{Pth}_{\boldsymbol{\mathcal{A}}^{(2)}}}_{s}\left(
\mathrm{ip}^{(2,[1])\sharp}_{s}\left(
\left[
\mathrm{CH}^{(1)}_{s}\left(
\mathfrak{P}
\right)\right]_{s}\right)
\right)\right)
\tag{1}
\\&=
\mathrm{CH}^{(2)}_{s}\left(
\mathrm{ip}^{(2,[1])\sharp}_{s}\left(
\mathrm{sc}^{0[\mathbf{PT}_{\boldsymbol{\mathcal{A}}}]}_{s}\left(
\mathrm{sc}^{([1],2)}_{s}\left(
\mathrm{ip}^{(2,[1])\sharp}_{s}\left(
\left[
\mathrm{CH}^{(1)}_{s}\left(
\mathfrak{P}
\right)\right]_{s}\right)
\right)
\right)\right)\right)
\tag{2}
\\&=
\mathrm{CH}^{(2)}_{s}\left(
\mathrm{ip}^{(2,[1])\sharp}_{s}\left(
\mathrm{sc}^{0[\mathbf{PT}_{\boldsymbol{\mathcal{A}}}]}_{s}\left(
\left[
\mathrm{CH}^{(1)}_{s}\left(
\mathfrak{P}
\right)\right]_{s}\right)
\right)
\right)
\tag{3}
\\&=
\eta^{(2,1)\sharp}_{s}\left(
\mathrm{CH}^{(1)\mathrm{m}}_{s}\left(
\mathrm{ip}^{([1],X)@}_{s}\left(
\mathrm{sc}^{0[\mathbf{PT}_{\boldsymbol{\mathcal{A}}}]}_{s}\left(
\left[
\mathrm{CH}^{(1)}_{s}\left(
\mathfrak{P}
\right)\right]_{s}\right)
\right)
\right)\right)
\tag{4}
\\&=
\eta^{(2,1)\sharp}_{s}\left(
\mathrm{CH}^{(1)}_{s}\left(
\mathrm{ip}^{(1,X)@}_{s}\left(
\mathrm{sc}^{0\mathbf{PT}_{\boldsymbol{\mathcal{A}}}}_{s}\left(
\mathrm{CH}^{(1)}_{s}\left(
\mathfrak{P}
\right)\right)
\right)
\right)\right)
\tag{5}
\\&=
\eta^{(2,1)\sharp}_{s}\left(
\mathrm{CH}^{(1)}_{s}\left(
\mathrm{sc}^{0\mathbf{Pth}_{\boldsymbol{\mathcal{A}}}}_{s}\left(
\mathrm{ip}^{(1,X)@}_{s}\left(
\mathrm{CH}^{(1)}_{s}\left(
\mathfrak{P}
\right)\right)
\right)
\right)\right)
\tag{6}
\\&=
\eta^{(2,1)\sharp}_{s}\left(
\mathrm{CH}^{(1)}_{s}\left(
\mathrm{ip}^{(1,0)\sharp}_{s}\left(
\mathrm{sc}^{(0,1)}_{s}\left(
\mathrm{ip}^{(1,X)@}_{s}\left(
\mathrm{CH}^{(1)}_{s}\left(
\mathfrak{P}
\right)\right)\right)
\right)
\right)\right)
\tag{7}
\\&=
\eta^{(2,1)\sharp}_{s}\left(
\mathrm{CH}^{(1)}_{s}\left(
\mathrm{ip}^{(1,0)\sharp}_{s}\left(
\mathrm{sc}^{(0,1)}_{s}\left(
\mathfrak{P}
\right)\right)
\right)\right).
\tag{8}
\\&=
\eta^{(2,1)\sharp}_{s}\left(
\mathrm{CH}^{(1)}_{s}\left(
\mathrm{sc}^{0\mathbf{Pth}_{\boldsymbol{\mathcal{A}}}}_{s}\left(
\mathfrak{P}
\right)
\right)\right).
\tag{9}
\end{align*}
\end{flushleft}

In the just stated chain of equalities, the first equality recovers the description of the second-order path $\mathfrak{P}^{(2)}$; the second equality follows from Claim~\ref{CDPthCatAlgScZ}; the third equality follows from Proposition~\ref{PDBasicEq}; the fourth equality follows from Proposition~\ref{PDCHDUId}; the fifth equality unravels the description of the mappings $\mathrm{CH}^{(1)\mathrm{m}}$, $\mathrm{ip}^{([1],X)@}$ and $\mathrm{sc}^{0[\mathbf{PT}_{\boldsymbol{\mathcal{A}}}]}_{s}$ according to Definition~\ref{DCHQuot}, Definition~\ref{DPTQIp} and Proposition~\ref{PPTQCatAlg}, respectively; the sixth equality follows from the fact that $\mathrm{ip}^{(1,X)@}$ is a $\Sigma^{\boldsymbol{\mathcal{A}}}$-homomorphism according to Definition~\ref{DIp}; the seventh equality unravels the description of the $0$-source operation symbol in the many-sorted partial $\Sigma^{\boldsymbol{\mathcal{A}}}$-algebra according to Proposition~\ref{PPthCatAlg}; the eighth equality follows from Proposition~\ref{PIpCH} and Lemma~\ref{LCH}; finally the last equality recovers the description of the $0$-source operation symbol in the many-sorted partial $\Sigma^{\boldsymbol{\mathcal{A}}}$-algebra according to Proposition~\ref{PPthCatAlg}.

Now, consider the right hand side of the pair in Equation~\ref{EDThetaDUI}. The following chain of equalities holds
\begin{flushleft}
$\mathrm{sc}^{0\mathbf{T}_{\Sigma^{\boldsymbol{\mathcal{A}}^{(2)}}}(X)}_{s}\left(
\mathrm{CH}^{(2)}_{s}\left(
\mathfrak{P}^{(2)}
\right)\right)$
\allowdisplaybreaks
\begin{align*}
&=\mathrm{sc}^{0\mathbf{T}_{\Sigma^{\boldsymbol{\mathcal{A}}^{(2)}}}(X)}_{s}\left(
\mathrm{CH}^{(2)}_{s}\left(
\mathrm{ip}^{(2,[1])\sharp}_{s}\left(
\left[
\mathrm{CH}^{(1)}_{s}\left(
\mathfrak{P}
\right)
\right]_{s}
\right)
\right)\right)
\tag{1}
\\&=
\mathrm{sc}^{0\mathbf{T}_{\Sigma^{\boldsymbol{\mathcal{A}}^{(2)}}}(X)}_{s}\left(
\eta^{(2,1)\sharp}_{s}\left(
\mathrm{CH}^{(1)\mathrm{m}}_{s}\left(
\mathrm{ip}^{([1],X)@}_{s}\left(
\left[
\mathrm{CH}^{(1)}_{s}\left(
\mathfrak{P}
\right)
\right]_{s}
\right)
\right)
\right)
\right)
\tag{2}
\\&=
\mathrm{sc}^{0\mathbf{T}_{\Sigma^{\boldsymbol{\mathcal{A}}^{(2)}}}(X)}_{s}\left(
\eta^{(2,1)\sharp}_{s}\left(
\mathrm{CH}^{(1)}_{s}\left(
\mathrm{ip}^{(1,X)@}_{s}\left(
\mathrm{CH}^{(1)}_{s}\left(
\mathfrak{P}
\right)
\right)
\right)
\right)
\right)
\tag{3}
\\&=
\mathrm{sc}^{0\mathbf{T}_{\Sigma^{\boldsymbol{\mathcal{A}}^{(2)}}}(X)}_{s}\left(
\eta^{(2,1)\sharp}_{s}\left(
\mathrm{CH}^{(1)}_{s}\left(
\mathfrak{P}
\right)
\right)
\right)
\tag{4}
\\&=
\eta^{(2,1)\sharp}_{s}\left(
\mathrm{sc}^{0\mathbf{T}_{\Sigma^{\boldsymbol{\mathcal{A}}}}(X)}_{s}\left(
\mathrm{CH}^{(1)}_{s}\left(
\mathfrak{P}
\right)
\right)
\right).
\tag{5}
\end{align*}
\end{flushleft}

In the just stated chain of equalities, the first equality recovers the description of the second-order path $\mathfrak{P}^{(2)}$; the second equality follows from Proposition~\ref{PDCHDUId}; the third equality unravels the description of the mappings $\mathrm{CH}^{(1)\mathrm{m}}$, $\mathrm{ip}^{([1],X)@}$ and $\mathrm{sc}^{0[\mathbf{PT}_{\boldsymbol{\mathcal{A}}}]}_{s}$ according to Definition~\ref{DCHQuot}, Definition~\ref{DPTQIp} and Proposition~\ref{PPTQCatAlg}, respectively; the fourth equality follows from Proposition~\ref{PIpCH}; finally, the last equality follows from Proposition~\ref{PDUEmb}.

This completes Case~(1).

\textsf{Case~(2)} Let $s$ be a sort in $S$ and $\mathfrak{P}$ a path in $\mathrm{Pth}_{\boldsymbol{\mathcal{A}},s}$, we want to prove that 
\[
\left(
\eta^{(2,1)\sharp}_{s}\left(
\mathrm{CH}^{(1)}_{s}\left(
\mathrm{tg}^{0\mathbf{Pth}_{\boldsymbol{\mathcal{A}}}}_{s}\left(
\mathfrak{P}
\right)\right)
\right)
,
\eta^{(2,1)\sharp}_{s}\left(
\mathrm{tg}^{0\mathbf{T}_{\Sigma^{\boldsymbol{\mathcal{A}}}}(X)}_{s}\left(
\mathrm{CH}^{(1)}_{s}\left(
\mathfrak{P}
\right)\right)
\right)
\right)
\in\Theta^{(2)}_{s}.
\]

For this, consider the second-order path $\mathfrak{P}^{(2)}=\mathrm{ip}^{(2,[1])\sharp}_{s}([\mathrm{CH}^{(1)}_{s}(\mathfrak{P})]_{s})$, according to Definition~\ref{DDTheta}, the following pair is in $\Theta^{(2)}_{s}$, 
\begin{equation}
\left(\mathrm{CH}^{(2)}_{s}\left(
\mathrm{tg}^{0\mathbf{Pth}_{\boldsymbol{\mathcal{A}}^{(2)}}}_{s}\left(
\mathfrak{P}^{(2)}
\right)\right),
\mathrm{tg}^{0\mathbf{T}_{\Sigma^{\boldsymbol{\mathcal{A}}^{(2)}}}(X)}_{s}\left(
\mathrm{CH}^{(2)}_{s}\left(
\mathfrak{P}^{(2)}
\right)\right)
\right)\in\Theta^{(2)}_{s}.
\tag{E2}\label{EDThetaDUII}
\end{equation}

This case follows by a similar argument to the one for the first case.

This completes Case~2.

\textsf{Case~(3)} Let $s$ be a sort in $S$ and $\mathfrak{Q}$ and $\mathfrak{P}$ two paths in $\mathrm{Pth}_{\boldsymbol{\mathcal{A}},s}$ satisfying that $\mathrm{sc}^{(0,1)}_{s}(\mathfrak{Q})=\mathrm{tg}^{(0,1)}_{s}(\mathfrak{P})$, we want to prove that 
\[
\left(
\eta^{(2,1)\sharp}_{s}\left(
\mathrm{CH}^{(1)}_{s}\left(
\mathfrak{Q}
\circ^{0\mathbf{Pth}_{\boldsymbol{\mathcal{A}}}}_{s}
\mathfrak{P}
\right)
\right)
,
\eta^{(2,1)\sharp}_{s}\left(
\mathrm{CH}^{(1)}_{s}\left(
\mathfrak{Q}
\right)
\circ^{0\mathbf{T}_{\Sigma^{\boldsymbol{\mathcal{A}}}}(X)}_{s}
\mathrm{CH}^{(1)}_{s}\left(
\mathfrak{P}
\right)
\right)
\right)
\in\Theta^{(2)}_{s}.
\]

For this, consider the second-order paths $\mathfrak{Q}^{(2)}=\mathrm{ip}^{(2,[1])\sharp}_{s}([\mathrm{CH}^{(1)}_{s}(\mathfrak{Q})]_{s})$ and $\mathfrak{P}^{(2)}=\mathrm{ip}^{(2,[1])\sharp}_{s}([\mathrm{CH}^{(1)}_{s}(\mathfrak{P})]_{s})$. First of all, note that the following chain of equalities holds
\allowdisplaybreaks
\begin{align*}
\mathrm{sc}^{(0,2)}_{s}\left(\mathfrak{Q}^{(2)}\right)&=
\mathrm{sc}^{(0,2)}_{s}\left(
\mathrm{ip}^{(2,[1])\sharp}_{s}\left(
\left[
\mathrm{CH}^{(1)}_{s}\left(
\mathfrak{Q}
\right)
\right]_{s}
\right)
\right)
\tag{1}
\\&=
\mathrm{sc}^{(0,[1])}_{s}\left(
\mathrm{ip}^{([1],X)@}_{s}\left(
\mathrm{sc}^{([1],2)}_{s}\left(
\mathrm{ip}^{(2,[1])\sharp}_{s}\left(
\left[
\mathrm{CH}^{(1)}_{s}\left(
\mathfrak{Q}
\right)
\right]_{s}
\right)
\right)\right)
\right)
\tag{2}
\\&=
\mathrm{sc}^{(0,[1])}_{s}\left(
\mathrm{ip}^{([1],X)@}_{s}\left(
\left[
\mathrm{CH}^{(1)}_{s}\left(
\mathfrak{Q}
\right)
\right]_{s}
\right)
\right)
\tag{3}
\\&=
\mathrm{sc}^{(0,[1])}_{s}\left(
\left[
\mathrm{ip}^{(1,X)@}_{s}\left(
\mathrm{CH}^{(1)}_{s}\left(
\mathfrak{Q}
\right)
\right)
\right]_{s}
\right)
\tag{4}
\\&=
\mathrm{sc}^{(0,[1])}_{s}\left(
\left[
\mathfrak{Q}
\right]_{s}
\right)
\tag{5}
\\&=
\mathrm{sc}^{(0,1)}_{s}\left(
\mathfrak{Q}
\right)
\tag{6}
\\&=
\mathrm{tg}^{(0,1)}_{s}\left(
\mathfrak{P}
\right)
\tag{7}
\\&=
\mathrm{tg}^{(0,[1])}_{s}\left(
\left[
\mathfrak{P}
\right]_{s}
\right)
\tag{8}
\\&=
\mathrm{tg}^{(0,[1])}_{s}\left(
\left[
\mathrm{ip}^{(1,X)@}_{s}\left(
\mathrm{CH}^{(1)}_{s}\left(
\mathfrak{P}
\right)
\right)
\right]_{s}
\right)
\tag{9}
\\&=
\mathrm{tg}^{(0,[1])}_{s}\left(
\mathrm{ip}^{([1],X)@}_{s}\left(
\left[
\mathrm{CH}^{(1)}_{s}\left(
\mathfrak{P}
\right)
\right]_{s}
\right)
\right)
\tag{10}
\\&=
\mathrm{tg}^{(0,[1])}_{s}\left(
\mathrm{ip}^{([1],X)@}_{s}\left(
\mathrm{tg}^{([1],2)}_{s}\left(
\mathrm{ip}^{(2,[1])\sharp}_{s}\left(
\left[
\mathrm{CH}^{(1)}_{s}\left(
\mathfrak{P}
\right)
\right]_{s}
\right)
\right)\right)
\right)
\tag{11}
\\&=
\mathrm{tg}^{(0,2)}_{s}\left(
\mathrm{ip}^{(2,[1])\sharp}_{s}\left(
\left[
\mathrm{CH}^{(1)}_{s}\left(
\mathfrak{P}
\right)
\right]_{s}
\right)
\right)
\tag{12}
\\&=
\mathrm{tg}^{(0,2)}_{s}\left(\mathfrak{P}^{(2)}\right).
\tag{13}
\end{align*}

In the just stated chain of equalities, the first equality recovers the description of the second-order path $\mathfrak{Q}^{(2)}$; the second equality unravels the description of the $(0,2)$-source mapping according to Definition~\ref{DDScTgZ}; the third equality follows from Proposition~\ref{PDBasicEq}; the fourth equality unravels the description of the mapping $\mathrm{ip}^{([1],X)@}$, according to Definition~\ref{DPTQIp}; the fifth equality follows from Proposition~\ref{PIpCH}; the sixth equality unravels the description of the mapping $\mathrm{sc}^{(0,[1])}$, according to Definition~\ref{DCHUZ}; the seventh equality follows by assumption; the eight equality unravels the description of the mapping $\mathrm{tg}^{(0,[1])}$, according to Definition~\ref{DCHUZ}; the ninth equality follows from Proposition~\ref{PIpCH}; the tenth equality unravels the description of the mapping $\mathrm{ip}^{([1],X)@}$, according to Definition~\ref{DPTQIp}; the eleventh equality follows from Proposition~\ref{PDBasicEq}; the twelfth equality recovers the description of the $(0,2)$-target mapping according to Definition~\ref{DDScTgZ}; finally, the last equality recovers the description of the second-order path $\mathfrak{P}^{(2)}$.

According to Definition~\ref{DDTheta}, the following pair is in $\Theta^{(2)}_{s}$, 
\allowdisplaybreaks
\begin{multline}
\left(\mathrm{CH}^{(2)}_{s}\left(
\mathfrak{Q}^{(2)}
\circ^{0\mathbf{Pth}_{\boldsymbol{\mathcal{A}}^{(2)}}}_{s}
\mathfrak{P}^{(2)}
\right),
\right.
\\
\left.
\mathrm{CH}^{(2)}_{s}\left(
\mathfrak{Q}^{(2)}
\right)
\circ^{0\mathbf{T}_{\Sigma^{\boldsymbol{\mathcal{A}}^{(2)}}}(X)}_{s}
\mathrm{CH}^{(2)}_{s}\left(
\mathfrak{P}^{(2)}
\right)
\right)\in\Theta^{(2)}_{s}.
\tag{E3}\label{EDThetaDUIII}
\end{multline}

Consider the left hand side of the pair in Equation~\ref{EDThetaDUIII}. The following chain of equalities holds
\begin{flushleft}
$\mathrm{CH}^{(2)}_{s}\left(
\mathfrak{Q}^{(2)}
\circ^{0\mathbf{Pth}_{\boldsymbol{\mathcal{A}}^{(2)}}}_{s}
\mathfrak{P}^{(2)}\right)$
\allowdisplaybreaks
\begin{align*}
&=
\mathrm{CH}^{(2)}_{s}\left(
\mathrm{ip}^{(2,[1])\sharp}_{s}\left(
\left[
\mathrm{CH}^{(1)}_{s}\left(
\mathfrak{Q}
\right)
\right]_{s}
\right)
\circ^{0\mathbf{Pth}_{\boldsymbol{\mathcal{A}}^{(2)}}}_{s}
\mathrm{ip}^{(2,[1])\sharp}_{s}\left(
\left[
\mathrm{CH}^{(1)}_{s}\left(
\mathfrak{P}
\right)
\right]_{s}
\right)
\right)
\tag{1}
\\&=
\mathrm{CH}^{(2)}_{s}\left(
\mathrm{ip}^{(2,[1])\sharp}_{s}\left(
\left[
\mathrm{CH}^{(1)}_{s}\left(
\mathfrak{Q}
\right)
\circ^{0\mathbf{PT}_{\boldsymbol{\mathcal{A}}}}_{s}
\mathrm{CH}^{(1)}_{s}\left(
\mathfrak{P}
\right)
\right]_{s}
\right)
\right)
\tag{2}
\\&=
\eta^{(2,1)\sharp}_{s}\left(
\mathrm{CH}^{(1)\mathrm{m}}_{s}\left(
\mathrm{ip}^{([1],X)@}_{s}\left(
\left[
\mathrm{CH}^{(1)}_{s}\left(
\mathfrak{Q}
\right)
\circ^{0\mathbf{PT}_{\boldsymbol{\mathcal{A}}}}_{s}
\mathrm{CH}^{(1)}_{s}\left(
\mathfrak{P}
\right)
\right]_{s}
\right)
\right)
\right)
\tag{3}
\\&=
\eta^{(2,1)\sharp}_{s}\left(
\mathrm{CH}^{(1)}_{s}\left(
\mathrm{ip}^{(1,X)@}_{s}\left(
\mathrm{CH}^{(1)}_{s}\left(
\mathfrak{Q}
\right)
\circ^{0\mathbf{PT}_{\boldsymbol{\mathcal{A}}}}_{s}
\mathrm{CH}^{(1)}_{s}\left(
\mathfrak{P}
\right)
\right)
\right)
\right)
\tag{4}
\\&=
\eta^{(2,1)\sharp}_{s}\left(
\mathrm{CH}^{(1)}_{s}\left(
\mathfrak{Q}
\circ^{0\mathbf{Pth}_{\boldsymbol{\mathcal{A}}}}_{s}
\mathfrak{P}
\right)
\right).
\tag{5}
\end{align*}
\end{flushleft}

In the just stated chain of equalities, the first equality recovers the description of the second-order paths $\mathfrak{Q}^{(2)}$ and  $\mathfrak{P}^{(2)}$; the second equality follows from Proposition~\ref{PDUIp}; the third equality follows from Proposition~\ref{PDCHDUId}; the fourth equality unravels the description of the mappings $\mathrm{CH}^{(1)\mathrm{m}}$ and $\mathrm{ip}^{([1],X)@}$ according to Definition~\ref{DCHQuot} and Definition~\ref{DPTQIp}, respectively; finally, the last equality follows from Lemma~\ref{LThetaComp}.

Now, consider the right hand side of the pair in Equation~\ref{EDThetaDUIII}. The following chain of equalities holds
\begin{flushleft}
$\mathrm{CH}^{(2)}_{s}\left(
\mathfrak{Q}^{(2)}
\right)
\circ^{0\mathbf{T}_{\Sigma^{\boldsymbol{\mathcal{A}}^{(2)}}}(X)}_{s}
\mathrm{CH}^{(2)}_{s}\left(
\mathfrak{P}^{(2)}
\right)$
\allowdisplaybreaks
\begin{align*}
\qquad&=
\mathrm{CH}^{(2)}_{s}\left(
\mathrm{ip}^{(2,[1])\sharp}_{s}\left(
\left[
\mathrm{CH}^{(1)}_{s}\left(
\mathfrak{Q}
\right)
\right]_{s}
\right)
\right)
\circ^{0\mathbf{T}_{\Sigma^{\boldsymbol{\mathcal{A}}^{(2)}}}(X)}_{s}
\\&\qquad\qquad\qquad\qquad\qquad\qquad\qquad
\mathrm{CH}^{(2)}_{s}\left(
\mathrm{ip}^{(2,[1])\sharp}_{s}\left(
\left[
\mathrm{CH}^{(1)}_{s}\left(
\mathfrak{P}
\right)
\right]_{s}
\right)
\right)
\tag{1}
\\&=
\eta^{(2,1)\sharp}_{s}\left(
\mathrm{CH}^{(1)\mathrm{m}}_{s}\left(
\mathrm{ip}^{([1],X)@}_{s}\left(
\left[
\mathrm{CH}^{(1)}_{s}\left(
\mathfrak{Q}
\right)
\right]_{s}
\right)\right)\right)
\circ^{0\mathbf{T}_{\Sigma^{\boldsymbol{\mathcal{A}}^{(2)}}}(X)}_{s}
\\&\qquad\qquad\qquad\qquad\qquad
\eta^{(2,1)\sharp}_{s}\left(
\mathrm{CH}^{(1)\mathrm{m}}_{s}\left(
\mathrm{ip}^{([1],X)@}_{s}\left(
\left[
\mathrm{CH}^{(1)}_{s}\left(
\mathfrak{P}
\right)
\right]_{s}
\right)\right)\right)
\tag{2}
\\&=
\eta^{(2,1)\sharp}_{s}\left(
\mathrm{CH}^{(1)}_{s}\left(
\mathrm{ip}^{(1,X)@}_{s}\left(
\mathrm{CH}^{(1)}_{s}\left(
\mathfrak{Q}
\right)
\right)\right)\right)
\circ^{0\mathbf{T}_{\Sigma^{\boldsymbol{\mathcal{A}}^{(2)}}}(X)}_{s}
\\&\qquad\qquad\qquad\qquad\qquad\qquad\qquad
\eta^{(2,1)\sharp}_{s}\left(
\mathrm{CH}^{(1)}_{s}\left(
\mathrm{ip}^{(1,X)@}_{s}\left(
\mathrm{CH}^{(1)}_{s}\left(
\mathfrak{P}
\right)
\right)\right)\right)
\tag{3}
\\&=
\eta^{(2,1)\sharp}_{s}\left(
\mathrm{CH}^{(1)}_{s}\left(
\mathfrak{Q}
\right)\right)
\circ^{0\mathbf{T}_{\Sigma^{\boldsymbol{\mathcal{A}}^{(2)}}}(X)}_{s}
\eta^{(2,1)\sharp}_{s}\left(
\mathrm{CH}^{(1)}_{s}\left(
\mathfrak{P}
\right)\right)
\tag{4}
\\&=
\eta^{(2,1)\sharp}_{s}\left(
\mathrm{CH}^{(1)}_{s}\left(
\mathfrak{Q}
\right)
\circ^{0\mathbf{T}_{\Sigma^{\boldsymbol{\mathcal{A}}}}(X)}_{s}
\mathrm{CH}^{(1)}_{s}\left(
\mathfrak{P}
\right)\right).
\tag{5}
\end{align*}
\end{flushleft}

In the just stated chain of equalities, the first equality recovers the description of the second-order paths $\mathfrak{Q}^{(2)}$ and $\mathfrak{P}^{(2)}$; the second equality follows from Proposition~\ref{PDCHDUId}; the third equality unravels the description of the mappings $\mathrm{CH}^{(1)\mathrm{m}}$ and $\mathrm{ip}^{([1],X)@}$,  according to Definition~\ref{DCHQuot} and Definition~\ref{DPTQIp}, respectively; the fourth equality follows from Proposition~\ref{PIpCH}; finally, the last equality follows from Proposition~\ref{PDUEmb}.

This completes Case~(3).

This completes the proof.
\end{proof}

We will continue studying the relation $\Theta^{(2)}$.
Let us recall that, for every sort $s\in S$, the term on the left defining $\Theta^{(2)}_{s}$ is always a term in $\mathrm{CH}^{(2)}_{s}[\mathrm{Pth}_{\boldsymbol{\mathcal{A}}^{(2)},s}]$, whilst the term on the right is not necessarily a term in such a subset of $\mathrm{T}_{\Sigma^{\boldsymbol{\mathcal{A}}^{(2)}}}(X)_{s}$.

However, we will show below that, for every sort $s\in S$, $\Theta^{(2)}_{s}$ relates pairs of terms in $\mathrm{T}_{\Sigma^{\boldsymbol{\mathcal{A}}^{(2)}}}(X)_{s}$ that, when mapped to $\mathrm{F}_{\Sigma^{\boldsymbol{\mathcal{A}}^{(2)}}}(\mathbf{Pth}_{\boldsymbol{\mathcal{A}}^{(2)}})_{s}$ under the action of $\mathrm{ip}^{(2,X)@}_{s}$ retrieve a pair of second-order paths that are in the same $\mathrm{Ker}(\mathrm{CH}^{(2)})_{s}$-class.

This fact already holds for the terms on the left of the pairs defining $\Theta^{(2)}_{s}$. In fact this is a direct consequence of Proposition~\ref{PDIpDCH}. Therefore, we will focus our attention exclusively on the terms appearing on the right of the pairs  defining $\Theta^{(2)}_{s}$.

\begin{lemma}\label{LDThetaScZ} 
Let $s$ be a sort in $S$ and $\mathfrak{P}^{(2)}$ a second-order path in $\mathrm{Pth}_{\boldsymbol{\mathcal{A}}^{(2)},s}$. Then
\begin{itemize}
\item[(i)] $\mathrm{ip}^{(2,X)@}_{s}(
\mathrm{sc}^{0\mathbf{T}_{\Sigma^{\boldsymbol{\mathcal{A}}^{(2)}}}(X)}_{s}(
\mathrm{CH}^{(2)}_{s}(
\mathfrak{P}^{(2)}
)))
$ is a second-order path in $\mathrm{Pth}_{\boldsymbol{\mathcal{A}}^{(2)},s}$;
\item[(ii)] $\mathrm{CH}^{(2)}_{s}(
\mathrm{ip}^{(2,X)@}_{s}(
\mathrm{sc}^{0\mathbf{T}_{\Sigma^{\boldsymbol{\mathcal{A}}^{(2)}}}(X)}_{s}(
\mathrm{CH}^{(2)}_{s}(
\mathfrak{P}^{(2)}
))))
=
\mathrm{CH}^{(2)}_{s}(
\mathrm{sc}^{0\mathbf{Pth}_{\boldsymbol{\mathcal{A}}^{(2)}}}_{s}(
\mathfrak{P}^{(2)}
))
$.
\end{itemize}
\end{lemma}
\begin{proof}
The following chain of equalities holds
\begin{flushleft}
$\mathrm{ip}^{(2,X)@}_{s}\left(
\mathrm{sc}^{0\mathbf{T}_{\Sigma^{\boldsymbol{\mathcal{A}}^{(2)}}}(X)}_{s}\left(
\mathrm{CH}^{(2)}_{s}\left(
\mathfrak{P}^{(2)}
\right)\right)\right)$
\allowdisplaybreaks
\begin{align*}
\quad
&=
\mathrm{sc}^{0\mathbf{F}_{\Sigma^{\boldsymbol{\mathcal{A}}^{(2)}}}
(\mathbf{Pth}_{\boldsymbol{\mathcal{A}}^{(2)}}
)}_{s}\left(
\mathrm{ip}^{(2,X)@}_{s}\left(
\mathrm{CH}^{(2)}_{s}\left(
\mathfrak{P}^{(2)}
\right)\right)\right)
\tag{1}
\\&=
\mathrm{sc}^{0\mathbf{Pth}_{\boldsymbol{\mathcal{A}}^{(2)}}}_{s}\left(
\mathrm{ip}^{(2,X)@}_{s}\left(
\mathrm{CH}^{(2)}_{s}\left(
\mathfrak{P}^{(2)}
\right)\right)\right)
\tag{2}
\\&=
\mathrm{ip}^{(2,0)\sharp}_{s}\left(
\mathrm{sc}^{(0,2)}_{s}\left(
\mathrm{ip}^{(2,X)@}_{s}\left(
\mathrm{CH}^{(2)}_{s}\left(
\mathfrak{P}^{(2)}
\right)\right)\right)\right)
\tag{3}
\\&=
\mathrm{ip}^{(2,0)\sharp}_{s}\left(
\mathrm{sc}^{(0,2)}_{s}\left(
\mathfrak{P}^{(2)}
\right)\right)
\tag{4}
\\&=
\mathrm{sc}^{0\mathbf{Pth}_{\boldsymbol{\mathcal{A}}^{(2)}}}_{s}\left(
\mathfrak{P}^{(2)}
\right).
\tag{5}
\end{align*}
\end{flushleft}

The first equality holds since $\mathrm{ip}^{(2,X)@}$ is a $\Sigma^{\boldsymbol{\mathcal{A}}^{(2)}}$-homomorphism in virtue of Definition~\ref{DDIp}; for the second equality, we note that, by Proposition~\ref{PDIpDCH} the element $\mathrm{ip}^{(2,X)@}_{s}(\mathrm{CH}^{(2)}_{s}(\mathfrak{P}^{(2)}))$ is a second-order path in $[\mathfrak{P}^{(2)}]^{}_{s}$, thus the interpretation of the operation symbol $\mathrm{sc}^{0}_{s}$ in the $\Sigma^{\boldsymbol{\mathcal{A}}^{(2)}}$-algebra $\mathbf{F}_{\Sigma^{\boldsymbol{\mathcal{A}}^{(2)}}}(\mathbf{Pth}_{\boldsymbol{\mathcal{A}}^{(2)}})$ is given by the respective interpretation of the operation symbol $\mathrm{sc}^{0}_{s}$ in the $\Sigma^{\boldsymbol{\mathcal{A}}^{(2)}}$-algebra $\mathbf{Pth}_{\boldsymbol{\mathcal{A}}^{(2)}}$; the third equality unravels the interpretation of the operation symbol $\mathrm{sc}^{0}_{s}$ in the $\Sigma^{\boldsymbol{\mathcal{A}}^{(2)}}$-algebra $\mathbf{Pth}_{\boldsymbol{\mathcal{A}}^{(2)}}$, as defined in Proposition~\ref{PDPthCatAlg}; the fourth equality follows from the fact that, by Proposition~\ref{PDIpDCH}, $\mathrm{ip}^{(2,X)@}_{s}(\mathrm{CH}^{(2)}_{s}(\mathfrak{P}^{(2)}))$ is a second-order path in $[\mathfrak{P}^{(2)}]^{}_{s}$, thus, by Corollary~\ref{CDCH} the second-order paths $\mathrm{ip}^{(2,X)@}_{s}(\mathrm{CH}^{(2)}_{s}(\mathfrak{P}^{(2)}))$ and $\mathfrak{P}^{(2)}$ have the same $(0,2)$-source; finally, the last equality unravels the interpretation of the operation symbol $\mathrm{sc}^{0}_{s}$ in the $\Sigma^{\boldsymbol{\mathcal{A}}^{(2)}}$-algebra $\mathbf{Pth}_{\boldsymbol{\mathcal{A}}^{(2)}}$. 

We note that, for a sort $s\in S$ and for a second-order path $\mathfrak{P}^{(2)}$ in $\mathrm{Pth}_{\boldsymbol{\mathcal{A}}^{(2)},s}$, it is always the case that $\mathrm{sc}^{0\mathbf{Pth}_{\boldsymbol{\mathcal{A}}^{(2)}}}_{s}(\mathfrak{P}^{(2)})$ is a second-order path in $\mathrm{Pth}_{\boldsymbol{\mathcal{A}}^{(2)},s}$.

From the just given equality, it immediately follows that 
$$
\mathrm{CH}^{(2)}_{s}\left(
\mathrm{ip}^{(2,X)@}_{s}\left(
\mathrm{sc}^{0\mathbf{T}_{\Sigma^{\boldsymbol{\mathcal{A}}^{(2)}}}(X)}_{s}\left(
\mathrm{CH}^{(2)}_{s}\left(
\mathfrak{P}^{(2)}
\right)\right)\right)\right)
=
\mathrm{CH}^{(2)}_{s}\left(
\mathrm{sc}^{0\mathbf{Pth}_{\boldsymbol{\mathcal{A}}^{(2)}}}_{s}\left(
\mathfrak{P}^{(2)}
\right)\right).
$$

This completes the proof.
\end{proof}

\begin{lemma}\label{LDThetaTgZ} 
Let $s$ be a sort in $S$ and $\mathfrak{P}^{(2)}$ a second-order path in $\mathrm{Pth}_{\boldsymbol{\mathcal{A}}^{(2)},s}$. Then
\begin{itemize}
\item[(i)] $\mathrm{ip}^{(2,X)@}_{s}(
\mathrm{tg}^{0\mathbf{T}_{\Sigma^{\boldsymbol{\mathcal{A}}^{(2)}}}(X)}_{s}(
\mathrm{CH}^{(2)}_{s}(
\mathfrak{P}^{(2)}
)))
$ is a second-order path in $\mathrm{Pth}_{\boldsymbol{\mathcal{A}}^{(2)},s}$;
\item[(ii)] $\mathrm{CH}^{(2)}_{s}(
\mathrm{ip}^{(2,X)@}_{s}(
\mathrm{tg}^{0\mathbf{T}_{\Sigma^{\boldsymbol{\mathcal{A}}^{(2)}}}(X)}_{s}(
\mathrm{CH}^{(2)}_{s}(
\mathfrak{P}^{(2)}
))))
=
\mathrm{CH}^{(2)}_{s}(
\mathrm{tg}^{0\mathbf{Pth}_{\boldsymbol{\mathcal{A}}^{(2)}}}_{s}(
\mathfrak{P}^{(2)}
))
$.
\end{itemize}
\end{lemma}

\begin{proof}
As it was the case in the previous lemma, it is easy to see that 
$$\mathrm{ip}^{(2,X)@}_{s}\left(
\mathrm{tg}^{0\mathbf{T}_{\Sigma^{\boldsymbol{\mathcal{A}}^{(2)}}}(X)}_{s}\left(
\mathrm{CH}^{(2)}_{s}\left(
\mathfrak{P}^{(2)}
\right)\right)\right)
=
\mathrm{tg}^{0\mathbf{Pth}_{\boldsymbol{\mathcal{A}}^{(2)}}}_{s}\left(
\mathfrak{P}^{(2)}
\right).$$

The proof is analogous to that presented in Lemma~\ref{LDThetaScZ}.
\end{proof}

\begin{lemma}\label{LDThetaCompZ} 
Let $s$ be a sort in $S$ and $\mathfrak{Q}^{(2)},\mathfrak{P}^{(2)}$ second-order paths in $\mathrm{Pth}_{\boldsymbol{\mathcal{A}}^{(2)},s}$ such that 
$\mathrm{sc}^{(0,2)}_{s}(\mathfrak{Q}^{(2)})=\mathrm{tg}^{(0,2)}_{s}(\mathfrak{P}^{(2)})$. Then
\begin{itemize}
\item[(i)] $\mathrm{ip}^{(2,X)@}_{s}(
\mathrm{CH}^{(2)}_{s}(
\mathfrak{Q}^{(2)}
)
\circ^{0\mathbf{T}_{\Sigma^{\boldsymbol{\mathcal{A}}^{(2)}}}(X)}_{s}
\mathrm{CH}^{(2)}_{s}(
\mathfrak{P}^{(2)}
))
$ is a second-order path in $\mathrm{Pth}_{\boldsymbol{\mathcal{A}}^{(2)},s}$;
\item[(ii)] 
\begin{flushleft}
$\mathrm{CH}^{(2)}_{s}(
\mathrm{ip}^{(2,X)@}_{s}(
\mathrm{CH}^{(2)}_{s}(
\mathfrak{Q}^{(2)}
)
\circ^{0\mathbf{T}_{\Sigma^{\boldsymbol{\mathcal{A}}^{(2)}}}(X)}_{s}
\mathrm{CH}^{(2)}_{s}(
\mathfrak{P}^{(2)}
)))$
\begin{flushright}
$=
\mathrm{CH}^{(2)}_{s}(
\mathfrak{Q}^{(2)}
\circ^{0\mathbf{Pth}_{\boldsymbol{\mathcal{A}}^{(2)}}}_{s}
\mathfrak{P}^{(2)}
).$
\end{flushright}
\end{flushleft}
\end{itemize}
\end{lemma}
\begin{proof}
The following chain of equalities holds
\begin{flushleft}
$\mathrm{ip}^{(2,X)@}_{s}\left(
\mathrm{CH}^{(2)}_{s}\left(
\mathfrak{Q}^{(2)}
\right)
\circ_{s}^{0\mathbf{T}_{\Sigma^{\boldsymbol{\mathcal{A}}^{(2)}}}(X)}
\mathrm{CH}^{(2)}_{s}\left(
\mathfrak{P}^{(2)}
\right)
\right)
$
\allowdisplaybreaks
\begin{align*}
\quad &=
\mathrm{ip}^{(2,X)@}_{s}\left(
\mathrm{CH}^{(2)}_{s}\left(
\mathfrak{Q}^{(2)}
\right)\right)
\circ^{0
\mathbf{F}_{\Sigma^{\boldsymbol{\mathcal{A}}^{(2)}}}\left(
\mathbf{Pth}_{\boldsymbol{\mathcal{A}}^{(2)}}\right)
}_{s}
\mathrm{ip}^{(2,X)@}_{s}\left(
\mathrm{CH}^{(2)}_{s}\left(
\mathfrak{P}^{(2)}
\right)\right)
\tag{1}
\\&=
\mathrm{ip}^{(2,X)@}_{s}\left(
\mathrm{CH}^{(2)}_{s}\left(
\mathfrak{Q}^{(2)}
\right)\right)
\circ^{0
\mathbf{Pth}_{\boldsymbol{\mathcal{A}}^{(2)}}
}_{s}
\mathrm{ip}^{(2,X)@}_{s}\left(
\mathrm{CH}^{(2)}_{s}\left(
\mathfrak{P}^{(2)}
\right)\right).
\tag{2}
\end{align*}
\end{flushleft}

The first equality holds since $\mathrm{ip}^{(2,X)@}$ is a $\Sigma^{\boldsymbol{\mathcal{A}}^{(2)}}$-homomorphism in virtue of Definition~\ref{DDIp}; the second equality follows from the fact that, by Proposition~\ref{PDIpDCH} the element $\mathrm{ip}^{(2,X)@}_{s}(
\mathrm{CH}^{(2)}_{s}(
\mathfrak{Q}^{(2)}
))$ is a second-order path in $[\mathfrak{Q}^{(2)}]^{}_{s}$ and the element $\mathrm{ip}^{(2,X)@}_{s}(
\mathrm{CH}^{(2)}_{s}(
\mathfrak{P}^{(2)}
))$ is a second-order path in $[\mathfrak{P}^{(2)}]^{}_{s}$. Moreover, by Corollary~\ref{CDCH}, the $(0,2)$-source of $\mathrm{ip}^{(2,X)@}_{s}(
\mathrm{CH}^{(2)}_{s}(
\mathfrak{Q}^{(2)}
))$ is the same as the $(0,2)$-source of $\mathfrak{Q}^{(2)}$ and the $(0,2)$-target of $\mathrm{ip}^{(2,X)@}_{s}(
\mathrm{CH}^{(2)}_{s}(
\mathfrak{P}^{(2)}
))$ is the same as the $(0,2)$-target of $\mathfrak{P}^{(2)}$. So, considering the foregoing, we have that 
\allowdisplaybreaks
\begin{align*}
\mathrm{sc}^{(0,2)}_{s}\left(
\mathrm{ip}^{(2,X)@}_{s}\left(
\mathrm{CH}^{(2)}_{s}\left(
\mathfrak{Q}^{(2)}
\right)\right)\right)
&=
\mathrm{sc}^{(0,2)}_{s}\left(
\mathfrak{Q}^{(2)}
\right)
\\&=
\mathrm{tg}^{(0,2)}_{s}\left(
\mathfrak{P}^{(2)}
\right)
\\&=
\mathrm{tg}^{(0,2)}_{s}\left(
\mathrm{ip}^{(2,X)@}_{s}\left(
\mathrm{CH}^{(2)}_{s}\left(
\mathfrak{P}^{(2)}
\right)\right)\right).
\end{align*}
Therefore, according to Proposition~\ref{PDPthCatAlg}, the $0$-composition of the second-order paths $\mathrm{ip}^{(2,X)@}_{s}(
\mathrm{CH}^{(2)}_{s}(
\mathfrak{Q}^{(2)}
))$ and $\mathrm{ip}^{(2,X)@}_{s}(
\mathrm{CH}^{(2)}_{s}(
\mathfrak{P}^{(2)}
))$ is granted to exist in the $\Sigma^{\boldsymbol{\mathcal{A}}^{(2)}}$-algebra $\mathbf{Pth}_{\boldsymbol{\mathcal{A}}^{(2)}}$, thus the interpretation of the operation symbol $\circ^{0}_{s}$ in the $\Sigma^{\boldsymbol{\mathcal{A}}^{(2)}}$-algebra $\mathbf{F}_{\Sigma^{\boldsymbol{\mathcal{A}}^{(2)}}}(\mathbf{Pth}_{\boldsymbol{\mathcal{A}}^{(2)}})$ is given by the respective interpretation of the operation symbol $\circ^{0}_{s}$ in the $\Sigma^{\boldsymbol{\mathcal{A}}^{(2)}}$-algebra $\mathbf{Pth}_{\boldsymbol{\mathcal{A}}^{(2)}}$. 

We note that, for the second-order paths $\mathrm{ip}^{(2,X)@}_{s}(
\mathrm{CH}^{(2)}_{s}(
\mathfrak{Q}^{(2)}
))$  and $\mathrm{ip}^{(2,X)@}_{s}(
\mathrm{CH}^{(2)}_{s}(
\mathfrak{P}^{(2)}
))$, since its $0$-composition is granted to exist, then its $0$-composition is a second-order path in $\mathrm{Pth}_{\boldsymbol{\mathcal{A}}^{(2)},s}$.

Now, regarding the second item in the statement of this lemma, we will consider two cases; either (1) both $\mathfrak{Q}^{(2)}$ and $\mathfrak{P}^{(2)}$ are $(2,[1])$-identity second-order paths or (2) either $\mathfrak{Q}^{(2)}$ or $\mathfrak{P}^{(2)}$ is not a $(2,[1])$-identity second-order paths.

For the first case, the following chain of equalities holds
\begin{flushleft}
$\mathrm{ip}^{(2,X)@}_{s}\left(
\mathrm{CH}^{(2)}_{s}\left(
\mathfrak{Q}^{(2)}
\right)
\circ_{s}^{0\mathbf{T}_{\Sigma^{\boldsymbol{\mathcal{A}}^{(2)}}}(X)}
\mathrm{CH}^{(2)}_{s}\left(
\mathfrak{P}^{(2)}
\right)
\right)
$
\allowdisplaybreaks
\begin{align*}
\qquad&=
\mathrm{ip}^{(2,X)@}_{s}\left(
\mathrm{CH}^{(2)}_{s}\left(
\mathfrak{Q}^{(2)}
\right)\right)
\circ^{0
\mathbf{Pth}_{\boldsymbol{\mathcal{A}}^{(2)}}
}_{s}
\mathrm{ip}^{(2,X)@}_{s}\left(
\mathrm{CH}^{(2)}_{s}\left(
\mathfrak{P}^{(2)}
\right)\right)
\tag{1}
\\&=
\mathfrak{Q}^{(2)}
\circ^{0
\mathbf{Pth}_{\boldsymbol{\mathcal{A}}^{(2)}}
}_{s}
\mathfrak{P}^{(2)}.
\tag{2}
\end{align*}
\end{flushleft}

In the just stated chain of equalities, the first equality follows from the discussion above, whilst the second equality follows from Proposition~\ref{PDIpDCH} and Corollary~\ref{CDCHUId}.

Therefore, since the second-order paths $\mathrm{ip}^{(2,X)@}_{s}(
\mathrm{CH}^{(2)}_{s}(
\mathfrak{Q}^{(2)}
)
\circ_{s}^{0\mathbf{T}_{\Sigma^{\boldsymbol{\mathcal{A}}^{(2)}}}(X)}
\mathrm{CH}^{(2)}_{s}(
\mathfrak{P}^{(2)}
)
)$ and $\mathfrak{Q}^{(2)}
\circ^{0
\mathbf{Pth}_{\boldsymbol{\mathcal{A}}^{(2)}}
}_{s}
\mathfrak{P}^{(2)}$ are equal, they have the same image under the second-order Curry-Howard mapping.

This completes Case~(1).

Now, regarding Case~(2), the following chain of equalities holds
\begin{flushleft}
$\mathrm{CH}^{(2)}_{s}\left(
\mathrm{ip}^{(2,X)@}_{s}\left(
\mathrm{CH}^{(2)}_{s}\left(
\mathfrak{Q}^{(2)}
\right)
\circ^{0\mathbf{T}_{\Sigma^{\boldsymbol{\mathcal{A}}^{(2)}}}(X)}_{s}
\mathrm{CH}^{(2)}_{s}\left(
\mathfrak{P}^{(2)}
\right)\right)\right)$
\allowdisplaybreaks
\begin{align*}
\qquad&=
\mathrm{CH}^{(2)}_{s}\left(
\mathrm{ip}^{(2,X)@}_{s}\left(
\mathrm{CH}^{(2)}_{s}\left(
\mathfrak{Q}^{(2)}
\circ^{0\mathbf{Pth}_{\boldsymbol{\mathcal{A}}^{(2)}}}_{s}
\mathfrak{P}^{(2)}
\right)\right)\right)
\tag{1}
\\&=
\mathrm{CH}^{(2)}_{s}\left(
\mathfrak{Q}^{(2)}
\circ^{0\mathbf{Pth}_{\boldsymbol{\mathcal{A}}^{(2)}}}_{s}
\mathfrak{P}^{(2)}
\right).
\tag{2}
\end{align*}
\end{flushleft}

In the just stated chain of equalities, the first equality recovers the definition of the image of the second-order Curry-Howard mapping at $\mathfrak{Q}^{(2)}
\circ^{0\mathbf{Pth}_{\boldsymbol{\mathcal{A}}^{(2)}}}_{s}
\mathfrak{P}^{(2)}$ according to Definition~\ref{DDCH}; and the second-equality follows from Proposition~\ref{PDIpDCH}.

This completes Case~(2).

This completes the proof.
\end{proof}

\begin{lemma}\label{LDThetaScU} 
Let $s$ be a sort in $S$ and $\mathfrak{P}^{(2)}$ a second-order path in $\mathrm{Pth}_{\boldsymbol{\mathcal{A}}^{(2)},s}$. Then
\begin{itemize}
\item[(i)] $\mathrm{ip}^{(2,X)@}_{s}(
\mathrm{sc}^{1\mathbf{T}_{\Sigma^{\boldsymbol{\mathcal{A}}^{(2)}}}(X)}_{s}(
\mathrm{CH}^{(2)}_{s}(
\mathfrak{P}^{(2)}
)))
$ is a second-order path in $\mathrm{Pth}_{\boldsymbol{\mathcal{A}}^{(2)},s}$;
\item[(ii)] $\mathrm{CH}^{(2)}_{s}(
\mathrm{ip}^{(2,X)@}_{s}(
\mathrm{sc}^{1\mathbf{T}_{\Sigma^{\boldsymbol{\mathcal{A}}^{(2)}}}(X)}_{s}(
\mathrm{CH}^{(2)}_{s}(
\mathfrak{P}^{(2)}
))))
=
\mathrm{CH}^{(2)}_{s}(
\mathrm{sc}^{1\mathbf{Pth}_{\boldsymbol{\mathcal{A}}^{(2)}}}_{s}(
\mathfrak{P}^{(2)}
))
$.
\end{itemize}
\end{lemma}
\begin{proof}
The following chain of equalities holds
\begin{flushleft}
$\mathrm{ip}^{(2,X)@}_{s}\left(
\mathrm{sc}^{1\mathbf{T}_{\Sigma^{\boldsymbol{\mathcal{A}}^{(2)}}}(X)}_{s}\left(
\mathrm{CH}^{(2)}_{s}\left(
\mathfrak{P}^{(2)}
\right)\right)\right)$
\allowdisplaybreaks
\begin{align*}
\quad
&=
\mathrm{sc}^{1\mathbf{F}_{\Sigma^{\boldsymbol{\mathcal{A}}^{(2)}}}
(\mathbf{Pth}_{\boldsymbol{\mathcal{A}}^{(2)}}
)}_{s}\left(
\mathrm{ip}^{(2,X)@}_{s}\left(
\mathrm{CH}^{(2)}_{s}\left(
\mathfrak{P}^{(2)}
\right)\right)\right)
\tag{1}
\\&=
\mathrm{sc}^{1\mathbf{Pth}_{\boldsymbol{\mathcal{A}}^{(2)}}}_{s}\left(
\mathrm{ip}^{(2,X)@}_{s}\left(
\mathrm{CH}^{(2)}_{s}\left(
\mathfrak{P}^{(2)}
\right)\right)\right)
\tag{2}
\\&=
\mathrm{ip}^{(2,[1])\sharp}_{s}\left(
\mathrm{sc}^{([1],2)}_{s}\left(
\mathrm{ip}^{(2,X)@}_{s}\left(
\mathrm{CH}^{(2)}_{s}\left(
\mathfrak{P}^{(2)}
\right)\right)\right)\right)
\tag{3}
\\&=
\mathrm{ip}^{(2,[1])\sharp}_{s}\left(
\mathrm{sc}^{([1],2)}_{s}\left(
\mathfrak{P}^{(2)}
\right)\right)
\tag{4}
\\&=
\mathrm{sc}^{1\mathbf{Pth}_{\boldsymbol{\mathcal{A}}^{(2)}}}_{s}\left(
\mathfrak{P}^{(2)}
\right).
\tag{5}
\end{align*}
\end{flushleft}

The first equality holds since $\mathrm{ip}^{(2,X)@}$ is a $\Sigma^{\boldsymbol{\mathcal{A}}^{(2)}}$-homomorphism in virtue of Definition~\ref{DDIp}; for the second equality, we note that, by Proposition~\ref{PDIpDCH} the element $\mathrm{ip}^{(2,X)@}_{s}(\mathrm{CH}^{(2)}_{s}(\mathfrak{P}^{(2)}))$ is a second-order path in $[\mathfrak{P}^{(2)}]^{}_{s}$, thus the interpretation of the operation symbol $\mathrm{sc}^{1}_{s}$ in the $\Sigma^{\boldsymbol{\mathcal{A}}^{(2)}}$-algebra $\mathbf{F}_{\Sigma^{\boldsymbol{\mathcal{A}}^{(2)}}}(\mathbf{Pth}_{\boldsymbol{\mathcal{A}}^{(2)}})$ is given by the respective interpretation of the operation symbol $\mathrm{sc}^{1}_{s}$ in the $\Sigma^{\boldsymbol{\mathcal{A}}^{(2)}}$-algebra $\mathbf{Pth}_{\boldsymbol{\mathcal{A}}^{(2)}}$; the third equality unravels the interpretation of the operation symbol $\mathrm{sc}^{1}_{s}$ in the $\Sigma^{\boldsymbol{\mathcal{A}}^{(2)}}$-algebra $\mathbf{Pth}_{\boldsymbol{\mathcal{A}}^{(2)}}$, as defined in Proposition~\ref{PDPthDCatAlg}; the fourth equality follows from the fact that, by Proposition~\ref{PDIpDCH}, $\mathrm{ip}^{(2,X)@}_{s}(\mathrm{CH}^{(2)}_{s}(\mathfrak{P}^{(2)}))$ is a second-order path in $[\mathfrak{P}^{(2)}]^{}_{s}$, thus, by Lemma~\ref{LDCH} the second-order paths $\mathrm{ip}^{(2,X)@}_{s}(\mathrm{CH}^{(2)}_{s}(\mathfrak{P}^{(2)}))$ and $\mathfrak{P}^{(2)}$ have the same $([1],2)$-source; finally, the last equality unravels the interpretation of the operation symbol $\mathrm{sc}^{1}_{s}$ in the $\Sigma^{\boldsymbol{\mathcal{A}}^{(2)}}$-algebra $\mathbf{Pth}_{\boldsymbol{\mathcal{A}}^{(2)}}$. 

We note that, for a sort $s\in S$ and for a second-order path $\mathfrak{P}^{(2)}$ in $\mathrm{Pth}_{\boldsymbol{\mathcal{A}}^{(2)},s}$, it is always the case that $\mathrm{sc}^{1\mathbf{Pth}_{\boldsymbol{\mathcal{A}}^{(2)}}}_{s}(\mathfrak{P}^{(2)})$ is a second-order path in $\mathrm{Pth}_{\boldsymbol{\mathcal{A}}^{(2)},s}$.

From the just given equality, it immediately follows that 
$$
\mathrm{CH}^{(2)}_{s}\left(
\mathrm{ip}^{(2,X)@}_{s}\left(
\mathrm{sc}^{1\mathbf{T}_{\Sigma^{\boldsymbol{\mathcal{A}}^{(2)}}}(X)}_{s}\left(
\mathrm{CH}^{(2)}_{s}\left(
\mathfrak{P}^{(2)}
\right)\right)\right)\right)
=
\mathrm{CH}^{(2)}_{s}\left(
\mathrm{sc}^{1\mathbf{Pth}_{\boldsymbol{\mathcal{A}}^{(2)}}}_{s}\left(
\mathfrak{P}^{(2)}
\right)\right).
$$

This completes the proof.
\end{proof}

\begin{lemma}\label{LDThetaTgU} 
Let $s$ be a sort in $S$ and $\mathfrak{P}^{(2)}$ a second-order path in $\mathrm{Pth}_{\boldsymbol{\mathcal{A}}^{(2)},s}$. Then
\begin{itemize}
\item[(i)] $\mathrm{ip}^{(2,X)@}_{s}(
\mathrm{tg}^{1\mathbf{T}_{\Sigma^{\boldsymbol{\mathcal{A}}^{(2)}}}(X)}_{s}(
\mathrm{CH}^{(2)}_{s}(
\mathfrak{P}^{(2)}
)))
$ is a second-order path in $\mathrm{Pth}_{\boldsymbol{\mathcal{A}}^{(2)},s}$;
\item[(ii)] $\mathrm{CH}^{(2)}_{s}(
\mathrm{ip}^{(2,X)@}_{s}(
\mathrm{tg}^{1\mathbf{T}_{\Sigma^{\boldsymbol{\mathcal{A}}^{(2)}}}(X)}_{s}(
\mathrm{CH}^{(2)}_{s}(
\mathfrak{P}^{(2)}
))))
=
\mathrm{CH}^{(2)}_{s}(
\mathrm{tg}^{1\mathbf{Pth}_{\boldsymbol{\mathcal{A}}^{(2)}}}_{s}(
\mathfrak{P}^{(2)}
))
$.
\end{itemize}
\end{lemma}

\begin{proof}
As it was the case in the previous lemma, it is easy to see that 
$$\mathrm{ip}^{(2,X)@}_{s}\left(
\mathrm{tg}^{1\mathbf{T}_{\Sigma^{\boldsymbol{\mathcal{A}}^{(2)}}}(X)}_{s}\left(
\mathrm{CH}^{(2)}_{s}\left(
\mathfrak{P}^{(2)}
\right)\right)\right)
=
\mathrm{tg}^{1\mathbf{Pth}_{\boldsymbol{\mathcal{A}}^{(2)}}}_{s}\left(
\mathfrak{P}^{(2)}
\right).$$

The proof is analogous to that presented in Lemma~\ref{LDThetaScU}.
\end{proof}

\begin{lemma}\label{LDThetaCompU} 
Let $s$ be a sort in $S$ and $\mathfrak{Q}^{(2)},\mathfrak{P}^{(2)}$ second-order paths in $\mathrm{Pth}_{\boldsymbol{\mathcal{A}}^{(2)},s}$ such that 
$\mathrm{sc}^{([1],2)}_{s}(\mathfrak{Q}^{(2)})=\mathrm{tg}^{([1],2)}_{s}(\mathfrak{P}^{(2)})$. Then
\begin{itemize}
\item[(i)] $\mathrm{ip}^{(2,X)@}_{s}(
\mathrm{CH}^{(2)}_{s}(
\mathfrak{Q}^{(2)}
)
\circ^{1\mathbf{T}_{\Sigma^{\boldsymbol{\mathcal{A}}^{(2)}}}(X)}_{s}
\mathrm{CH}^{(2)}_{s}(
\mathfrak{P}^{(2)}
))
$ is a second-order path in $\mathrm{Pth}_{\boldsymbol{\mathcal{A}}^{(2)},s}$;
\item[(ii)] 
\begin{flushleft}
$\mathrm{CH}^{(2)}_{s}(
\mathrm{ip}^{(2,X)@}_{s}(
\mathrm{CH}^{(2)}_{s}(
\mathfrak{Q}^{(2)}
)
\circ^{1\mathbf{T}_{\Sigma^{\boldsymbol{\mathcal{A}}^{(2)}}}(X)}_{s}
\mathrm{CH}^{(2)}_{s}(
\mathfrak{P}^{(2)}
)))$
\begin{flushright}
$=
\mathrm{CH}^{(2)}_{s}(
\mathfrak{Q}^{(2)}
\circ^{1\mathbf{Pth}_{\boldsymbol{\mathcal{A}}^{(2)}}}_{s}
\mathfrak{P}^{(2)}
).$
\end{flushright}
\end{flushleft}
\end{itemize}
\end{lemma}

\begin{proof}
The following chain of equalities holds
\begin{flushleft}
$\mathrm{ip}^{(2,X)@}_{s}\left(
\mathrm{CH}^{(2)}_{s}\left(
\mathfrak{Q}^{(2)}
\right)
\circ_{s}^{1\mathbf{T}_{\Sigma^{\boldsymbol{\mathcal{A}}^{(2)}}}(X)}
\mathrm{CH}^{(2)}_{s}\left(
\mathfrak{P}^{(2)}
\right)
\right)
$
\allowdisplaybreaks
\begin{align*}
\quad &=
\mathrm{ip}^{(2,X)@}_{s}\left(
\mathrm{CH}^{(2)}_{s}\left(
\mathfrak{Q}^{(2)}
\right)\right)
\circ^{1
\mathbf{F}_{\Sigma^{\boldsymbol{\mathcal{A}}^{(2)}}}\left(
\mathbf{Pth}_{\boldsymbol{\mathcal{A}}^{(2)}}\right)
}_{s}
\mathrm{ip}^{(2,X)@}_{s}\left(
\mathrm{CH}^{(2)}_{s}\left(
\mathfrak{P}^{(2)}
\right)\right)
\tag{1}
\\&=
\mathrm{ip}^{(2,X)@}_{s}\left(
\mathrm{CH}^{(2)}_{s}\left(
\mathfrak{Q}^{(2)}
\right)\right)
\circ^{1
\mathbf{Pth}_{\boldsymbol{\mathcal{A}}^{(2)}}
}_{s}
\mathrm{ip}^{(2,X)@}_{s}\left(
\mathrm{CH}^{(2)}_{s}\left(
\mathfrak{P}^{(2)}
\right)\right).
\tag{2}
\end{align*}
\end{flushleft}

The first equality holds since $\mathrm{ip}^{(2,X)@}$ is a $\Sigma^{\boldsymbol{\mathcal{A}}^{(2)}}$-homomorphism in virtue of Definition~\ref{DDIp}; the second equality follows from the fact that, by Proposition~\ref{PDIpDCH} the element $\mathrm{ip}^{(2,X)@}_{s}(
\mathrm{CH}^{(2)}_{s}(
\mathfrak{Q}^{(2)}
))$ is a second-order path in $[\mathfrak{Q}^{(2)}]^{}_{s}$ and the element $\mathrm{ip}^{(2,X)@}_{s}(
\mathrm{CH}^{(2)}_{s}(
\mathfrak{P}^{(2)}
))$ is a second-order path in $[\mathfrak{P}^{(2)}]^{}_{s}$. Moreover, by Lemma~\ref{LDCH}, the $([1],2)$-source of $\mathrm{ip}^{(2,X)@}_{s}(
\mathrm{CH}^{(2)}_{s}(
\mathfrak{Q}^{(2)}
))$ is the same as the $([1],2)$-source of $\mathfrak{Q}^{(2)}$ and the $([1],2)$-target of $\mathrm{ip}^{(2,X)@}_{s}(
\mathrm{CH}^{(2)}_{s}(
\mathfrak{P}^{(2)}
))$ is the same as the $([1],2)$-target of $\mathfrak{P}^{(2)}$. So, considering the foregoing, we have that 
\allowdisplaybreaks
\begin{align*}
\mathrm{sc}^{([1],2)}_{s}\left(
\mathrm{ip}^{(2,X)@}_{s}\left(
\mathrm{CH}^{(2)}_{s}\left(
\mathfrak{Q}^{(2)}
\right)\right)\right)
&=
\mathrm{sc}^{([1],2)}_{s}\left(
\mathfrak{Q}^{(2)}
\right)
\\&=
\mathrm{tg}^{([1],2)}_{s}\left(
\mathfrak{P}^{(2)}
\right)
\\&=
\mathrm{tg}^{([1],2)}_{s}\left(
\mathrm{ip}^{(2,X)@}_{s}\left(
\mathrm{CH}^{(2)}_{s}\left(
\mathfrak{P}^{(2)}
\right)\right)\right).
\end{align*}
Therefore, according to Proposition~\ref{PDPthDCatAlg}, the $1$-composition of the second-order paths $\mathrm{ip}^{(2,X)@}_{s}(
\mathrm{CH}^{(2)}_{s}(
\mathfrak{Q}^{(2)}
))$ and $\mathrm{ip}^{(2,X)@}_{s}(
\mathrm{CH}^{(2)}_{s}(
\mathfrak{P}^{(2)}
))$ is granted to exist in the $\Sigma^{\boldsymbol{\mathcal{A}}^{(2)}}$-algebra $\mathbf{Pth}_{\boldsymbol{\mathcal{A}}^{(2)}}$, thus the interpretation of the operation symbol $\circ^{1}_{s}$ in the $\Sigma^{\boldsymbol{\mathcal{A}}^{(2)}}$-algebra $\mathbf{F}_{\Sigma^{\boldsymbol{\mathcal{A}}^{(2)}}}(\mathbf{Pth}_{\boldsymbol{\mathcal{A}}^{(2)}})$ is given by the respective interpretation of the operation symbol $\circ^{1}_{s}$ in the $\Sigma^{\boldsymbol{\mathcal{A}}^{(2)}}$-algebra $\mathbf{Pth}_{\boldsymbol{\mathcal{A}}^{(2)}}$. 

For the second-order paths $\mathrm{ip}^{(2,X)@}_{s}(
\mathrm{CH}^{(2)}_{s}(
\mathfrak{Q}^{(2)}
))$  and $\mathrm{ip}^{(2,X)@}_{s}(
\mathrm{CH}^{(2)}_{s}(
\mathfrak{P}^{(2)}
))$, since its $1$-composition is granted to exist, then its $1$-composition is a second-order path in $\mathrm{Pth}_{\boldsymbol{\mathcal{A}}^{(2)},s}$.

The second item in the statement of this lemma is done by Artinian induction on $(\coprod\mathrm{Pth}_{\boldsymbol{\mathcal{A}}^{(2)}}, \leq_{\mathbf{Pth}_{\boldsymbol{\mathcal{A}}^{(2)}}})$. But before considering the general case, we first consider the case in which one of the second-order paths is a $(2,[1])$-identity second-order path.

\begin{claim}\label{CDThetaCompU} Let $s$ be a sort in $S$ and $\mathfrak{Q}^{(2)},\mathfrak{P}^{(2)}$ second-order paths in $\mathrm{Pth}_{\boldsymbol{\mathcal{A}}^{(2)},s}$ such that $\mathrm{sc}^{([1],2)}_{s}(\mathfrak{Q}^{(2)})=\mathrm{tg}^{([1],2)}_{s}(\mathfrak{P}^{(2)})$. If either $\mathfrak{P}^{(2)}$ or $\mathfrak{Q}^{(2)}$ is a $(2,[1])$-identity second-order path, then
\begin{multline*}
\mathrm{CH}^{(2)}_{s}\left(
\mathrm{ip}^{(2,X)@}_{s}\left(
\mathrm{CH}^{(2)}_{s}\left(
\mathfrak{Q}^{(2)}
\right)
\circ^{1\mathbf{T}_{\Sigma^{\boldsymbol{\mathcal{A}}^{(2)}}}(X)}_{s}
\mathrm{CH}^{(2)}_{s}\left(
\mathfrak{P}^{(2)}
\right)\right)\right)
\\=
\mathrm{CH}^{(2)}_{s}\left(
\mathfrak{Q}^{(2)}
\circ_{s}^{1\mathbf{Pth}_{\boldsymbol{\mathcal{A}}^{(2)}}}
\mathfrak{P}^{(2)}
\right).
\end{multline*}
\end{claim}

Assume without loss of generality that $\mathfrak{P}^{(2)}$ is a $(2,[1])$-identity second-order path. Note that in this case the following chain of equalities holds
\begin{flushleft}
$\mathrm{CH}^{(2)}_{s}\left(
\mathrm{ip}^{(2,X)@}_{s}\left(
\mathrm{CH}^{(2)}_{s}\left(
\mathfrak{Q}^{(2)}
\right)
\circ^{1\mathbf{T}_{\Sigma^{\boldsymbol{\mathcal{A}}^{(2)}}}(X)}_{s}
\mathrm{CH}^{(2)}_{s}\left(
\mathfrak{P}^{(2)}
\right)\right)\right)
$
\allowdisplaybreaks
\begin{align*}
\quad &=
\mathrm{CH}^{(2)}_{s}\left(
\mathrm{ip}^{(2,X)@}_{s}\left(
\mathrm{CH}^{(2)}_{s}\left(
\mathfrak{Q}^{(2)}
\right)\right)
\circ^{1\mathbf{Pth}_{\boldsymbol{\mathcal{A}}^{(2)}}}_{s}
\mathrm{ip}^{(2,X)@}_{s}\left(
\mathrm{CH}^{(2)}_{s}\left(
\mathfrak{P}^{(2)}
\right)\right)\right)
\tag{1}
\\&=
\mathrm{CH}^{(2)}_{s}\left(
\mathrm{ip}^{(2,X)@}_{s}\left(
\mathrm{CH}^{(2)}_{s}\left(
\mathfrak{Q}^{(2)}
\right)\right)
\circ^{1\mathbf{Pth}_{\boldsymbol{\mathcal{A}}^{(2)}}}_{s}
\mathfrak{P}^{(2)}
\right)
\tag{2}
\\&=
\mathrm{CH}^{(2)}_{s}\left(
\mathrm{ip}^{(2,X)@}_{s}\left(
\mathrm{CH}^{(2)}_{s}\left(
\mathfrak{Q}^{(2)}
\right)\right)\right)
\tag{3}
\\&=
\mathrm{CH}^{(2)}_{s}\left(
\mathfrak{Q}^{(2)}
\right)
\tag{4}
\\&=
\mathrm{CH}^{(2)}_{s}\left(
\mathfrak{Q}^{(2)}
\circ^{1\mathbf{Pth}_{\boldsymbol{\mathcal{A}}^{(2)}}}_{s}
\mathfrak{P}^{(2)}
\right).
\tag{5}
\end{align*}
\end{flushleft}

The first equality is a direct consequence of the equation proven in the first item of this lemma; the second equality follows from Corollary~\ref{CDCHUId}, since $\mathfrak{P}^{(2)}$ is a $(2,[1])$-identity second-order path and by Definition~\ref{DDIp}, $\mathrm{ip}^{(2,X)@}_{s}(
\mathrm{CH}^{(2)}_{s}(
\mathfrak{P}^{(2)}
))$ is a second-order path in $[\mathfrak{P}^{(2)}]^{}_{s}$, therefore $\mathrm{ip}^{(2,X)@}_{s}(
\mathrm{CH}^{(2)}_{s}(
\mathfrak{P}^{(2)}
))=\mathfrak{P}^{(2)}$; the third equality follows from the fact that $\mathfrak{P}^{(2)}$ is a $(2,[1])$-identity second-order path and, by assumption $\mathrm{sc}^{([1],2)}_{s}(\mathfrak{Q}^{(2)})=\mathrm{tg}^{([1],2)}_{s}(\mathfrak{P}^{(2)})$, thus, by Proposition~\ref{PDUIp}, $\mathfrak{P}^{(2)}$ is a the $(2,[1])$-identity second-order path on $\mathrm{sc}^{([1],2)}_{s}(\mathfrak{Q}^{(2)})$. Moreover, by Proposition~\ref{PDIpDCH}, the element $\mathrm{ip}^{(2,X)@}_{s}(
\mathrm{CH}^{(2)}_{s}(
\mathfrak{Q}^{(2)}
))$ is a second-order path in $[\mathfrak{Q}^{(2)}]^{}_{s}$ and, by Lemma~\ref{LDCH}, the $([1],2)$-source of $\mathrm{ip}^{(2,X)@}_{s}(
\mathrm{CH}^{(2)}_{s}(
\mathfrak{Q}^{(2)}
))$ is the same as the $([1],2)$-source of $\mathfrak{Q}^{(2)}$. In the light of the above, we have that $$\mathrm{ip}^{(2,X)@}_{s}\left(
\mathrm{CH}^{(2)}_{s}\left(
\mathfrak{Q}^{(2)}
\right)\right)
\circ^{1\mathbf{Pth}_{\boldsymbol{\mathcal{A}}^{(2)}}}_{s}
\mathfrak{P}^{(2)}=
\mathrm{ip}^{(2,X)@}_{s}\left(
\mathrm{CH}^{(2)}_{s}\left(
\mathfrak{Q}^{(2)}
\right)\right);$$
the fourth equality follows from the fact that, by Proposition~\ref{PDIpDCH}, the element $\mathrm{ip}^{(2,X)@}_{s}(
\mathrm{CH}^{(2)}_{s}(
\mathfrak{Q}^{(2)}
))$ is a second-order path in $[\mathfrak{Q}^{(2)}]^{}_{s}$, thus its image under the second-order Curry-Howard mapping is the same as that of $\mathfrak{Q}^{(2)}$. Finally, the last equality follows from the fact that $\mathfrak{P}^{(2)}$ is the $(2,[1])$-identity second-order path on the $([1],2)$-source of $\mathfrak{Q}^{(2)}$ as we have mentioned before.

The same argument applies in case $\mathfrak{Q}^{(2)}$ is a $(2,[1])$-identity second-order path. In this case we will argue taking into consideration that $\mathfrak{Q}^{(2)}$ must be the $(2,[1])$-identity second-order path on the $([1],2)$-target of $\mathfrak{P}^{(2)}$.

This completes the proof of Claim~\ref{CDThetaCompU}.

We now prove the general case by Artinian induction on $(\coprod\mathrm{Pth}_{\boldsymbol{\mathcal{A}}^{(2)}}, \leq_{\mathbf{Pth}_{\boldsymbol{\mathcal{A}}^{(2)}}})$.

\textsf{Base step of the Artinian induction.}

Let $(\mathfrak{Q}^{(2)}\circ^{1\mathbf{Pth}_{\boldsymbol{\mathcal{A}}^{(2)}}}_{s}\mathfrak{P}^{(2)},s)$ be a minimal element in $(\coprod\mathrm{Pth}_{\boldsymbol{\mathcal{A}}^{(2)}},\leq_{\mathbf{Pth}_{\boldsymbol{\mathcal{A}}^{(2)}}})$. Then by Proposition~\ref{PDMinimal} the second-order path $\mathfrak{Q}^{(2)}\circ^{1\mathbf{Pth}_{\boldsymbol{\mathcal{A}}^{(2)}}}_{s}\mathfrak{P}^{(2)}$ is either a $(2,[1])$-identity second-order path or a   second-order echelon. In any case, either $\mathfrak{P}^{(2)}$ or $\mathfrak{Q}^{(2)}$ must be a $(2,[1])$-identity second-order path. The statement follows by Claim~\ref{CDThetaCompU}.

\textsf{Inductive step of the Artinian induction.}

Let $(\mathfrak{Q}^{(2)}\circ^{1\mathbf{Pth}_{\boldsymbol{\mathcal{A}}^{(2)}}}_{s}\mathfrak{P}^{(2)},s)$ be a non-minimal element in $(\coprod\mathrm{Pth}_{\boldsymbol{\mathcal{A}}^{(2)}},\leq_{\mathbf{Pth}_{\boldsymbol{\mathcal{A}}^{(2)}}})$. Let us suppose that, for every sort $t\in S$ and every second-order path $\mathfrak{Q}'^{(2)}\circ^{1\mathbf{Pth}_{\boldsymbol{\mathcal{A}}^{(2)}}}_{t}\mathfrak{P}'^{(2)}$ in $\mathrm{Pth}_{\boldsymbol{\mathcal{A}}^{(2)},t}$, if $(\mathfrak{Q}'^{(2)}\circ^{1\mathbf{Pth}_{\boldsymbol{\mathcal{A}}^{(2)}}}_{t}\mathfrak{P}'^{(2)},t)
<_{\mathbf{Pth}_{\boldsymbol{\mathcal{A}}^{(2)}}} (\mathfrak{Q}^{(2)}\circ^{1\mathbf{Pth}_{\boldsymbol{\mathcal{A}}^{(2)}}}_{s}\mathfrak{P}^{(2)},s)
$, then the statement holds for $\mathfrak{Q}'^{(2)}\circ^{1\mathbf{Pth}_{\boldsymbol{\mathcal{A}}^{(2)}}}_{t}\mathfrak{P}'^{(2)}$, i.e., the following equality holds
\begin{multline*}
\mathrm{CH}^{(2)}_{t}\left(
\mathrm{ip}^{(2,X)@}_{t}\left(
\mathrm{CH}^{(2)}_{t}\left(
\mathfrak{Q}'^{(2)}
\right)
\circ^{1\mathbf{T}_{\Sigma^{\boldsymbol{\mathcal{A}}^{(2)}}}(X)}_{t}
\mathrm{CH}^{(2)}_{t}
\left(\mathfrak{P}'^{(2)}
\right)\right)\right)
\\=
\mathrm{CH}^{(2)}_{t}\left(
\mathfrak{Q}'^{(2)}
\circ^{1\mathbf{Pth}_{\boldsymbol{\mathcal{A}}^{(2)}}}_{t}
\mathfrak{P}'^{(2)}
\right).
\end{multline*}

Since $(\mathfrak{Q}^{(2)}\circ^{1\mathbf{Pth}_{\boldsymbol{\mathcal{A}}^{(2)}}}_{s}\mathfrak{P}^{(2)},s)$ is a non-minimal element in $(\coprod\mathrm{Pth}_{\boldsymbol{\mathcal{A}}^{(2)}},\leq_{\mathbf{Pth}_{\boldsymbol{\mathcal{A}}^{(2)}}})$ and taking into account, in virtue of Claim~\ref{CDThetaCompU}, that we can assume that neither $\mathfrak{Q}^{(2)}$ nor $\mathfrak{P}^{(2)}$ are $(2,[1])$-identity second-order paths, then, by Lemma~\ref{LDOrdI}, $\mathfrak{Q}^{(2)}\circ^{1\mathbf{Pth}_{\boldsymbol{\mathcal{A}}^{(2)}}}_{s}\mathfrak{P}^{(2)}$ is either (1) a second-order path of length strictly greater than one containing at least one   second-order echelon or (2) an echelonless second-order path.

If~(1), then let $i\in\bb{\mathfrak{Q}^{(2)}\circ^{1\mathbf{Pth}_{\boldsymbol{\mathcal{A}}^{(2)}}}_{s}\mathfrak{P}^{(2)}}$ be the first index for which the one-step subpath $(\mathfrak{Q}^{(2)}\circ^{1\mathbf{Pth}_{\boldsymbol{\mathcal{A}}^{(2)}}}_{s}\mathfrak{P}^{(2)})^{i,i}$ of $\mathfrak{Q}^{(2)}\circ^{1\mathbf{Pth}_{\boldsymbol{\mathcal{A}}^{(2)}}}_{s}\mathfrak{P}^{(2)}$ is a   second-order echelon. We distinguish the following cases according to the nature of $i$; either~(1.1) $i=0$ or~(1.2) $i>0$.

If~(1.1), i.e., if $\mathfrak{Q}^{(2)}\circ^{1\mathbf{Pth}_{\boldsymbol{\mathcal{A}}^{(2)}}}_{s}\mathfrak{P}^{(2)}$ is a second-order path containing a   second-order echelon on its first step, since we are assuming that $\mathfrak{P}^{(2)}$ is not a $(2,[1])$-identity second-order path, then $\mathfrak{P}^{(2)}$ has a   second-order echelon on its first step. Then it could be the case that either (1.1.1) $\mathfrak{P}^{(2)}$ is a   second-order echelon or (1.1.2) $\mathfrak{P}^{(2)}$ is a second-order path of length strictly greater than one containing a   second-order echelon on its first step.

If~(1.1.1), i.e., if $\mathfrak{Q}^{(2)}\circ^{1\mathbf{Pth}_{\boldsymbol{\mathcal{A}}^{(2)}}}_{s}\mathfrak{P}^{(2)}$ is a second-order path containing a   second-order echelon on its first step and $\mathfrak{P}^{(2)}$ is a   second-order echelon then, according to Definition~\ref{DDCH}, the value of the second-order Curry-Howard mapping at $\mathfrak{Q}^{(2)}\circ^{1\mathbf{Pth}_{\boldsymbol{\mathcal{A}}^{(2)}}}_{s}\mathfrak{P}^{(2)}$ is given by
$$
\mathrm{CH}^{(2)}_{s}\left(
\mathfrak{Q}^{(2)}\circ^{1\mathbf{Pth}_{\boldsymbol{\mathcal{A}}^{(2)}}}_{s}\mathfrak{P}^{(2)}
\right)
=
\mathrm{CH}^{(2)}\left(
\mathfrak{Q}^{(2)}
\right)
\circ^{1\mathbf{T}_{\Sigma^{\boldsymbol{\mathcal{A}}^{(2)}}}(X)}_{s}
\mathrm{CH}^{(2)}_{s}\left(
\mathfrak{P}^{(2)}
\right).
$$

On the other hand, the following chain of equalities holds
\begin{flushleft}
$
\mathrm{CH}^{(2)}_{s}\left(
\mathrm{ip}^{(2,X)@}_{s}\left(
\mathrm{CH}^{(2)}_{s}\left(
\mathfrak{Q}^{(2)}
\right)
\circ^{1\mathbf{T}_{\Sigma^{\boldsymbol{\mathcal{A}}^{(2)}}}(X)}_{s}
\mathrm{CH}^{(2)}_{s}\left(
\mathfrak{P}^{(2)}
\right)\right)\right)
$
\allowdisplaybreaks
\begin{align*}
\quad&=
\mathrm{CH}^{(2)}_{s}\left(
\mathrm{ip}^{(2,X)@}_{s}\left(
\mathrm{CH}^{(2)}_{s}\left(
\mathfrak{Q}^{(2)}
\right)\right)
\circ^{1\mathbf{Pth}_{\boldsymbol{\mathcal{A}}^{(2)}}}_{s}
\mathrm{ip}^{(2,X)@}_{s}\left(
\mathrm{CH}^{(2)}_{s}\left(
\mathfrak{P}^{(2)}
\right)\right)\right)
\tag{1}
\\&=
\mathrm{CH}^{(2)}_{s}\left(
\mathrm{ip}^{(2,X)@}_{s}\left(
\mathrm{CH}^{(2)}_{s}\left(
\mathfrak{Q}^{(2)}
\right)\right)\right)
\circ^{1\mathbf{T}_{\Sigma^{\boldsymbol{\mathcal{A}}^{(2)}}}(X)}_{s}
\\&\qquad\qquad\qquad\qquad\qquad\qquad\qquad\qquad\qquad
\mathrm{CH}^{(2)}_{s}\left(
\mathrm{ip}^{(2,X)@}_{s}\left(
\mathrm{CH}^{(2)}_{s}\left(
\mathfrak{P}^{(2)}
\right)\right)\right)
\tag{2}
\\&=
\mathrm{CH}^{(2)}_{s}\left(
\mathfrak{Q}^{(2)}
\right)
\circ^{1\mathbf{T}_{\Sigma^{\boldsymbol{\mathcal{A}}^{(2)}}}(X)}_{s}
\mathrm{CH}^{(2)}_{s}\left(
\mathfrak{P}^{(2)}
\right)
\tag{3}
\\&=
\mathrm{CH}^{(2)}_{s}\left(
\mathfrak{Q}^{(2)}\circ^{1\mathbf{Pth}_{\boldsymbol{\mathcal{A}}^{(2)}}}_{s}\mathfrak{P}^{(2)}
\right).
\tag{4}
\end{align*}
\end{flushleft}

The first equality is a direct consequence of the equation proven in the first item of this lemma; the second equality follows from the fact that, by Proposition~\ref{PDIpDCH}, the element $\mathrm{ip}^{(2,X)@}_{s}(
\mathrm{CH}^{(2)}_{s}(
\mathfrak{P}^{(2)}
))$ is a second-order path in $[\mathfrak{P}^{(2)}]^{}_{s}$ therefore, since we are assuming that $\mathfrak{P}^{(2)}$ is a   second-order echelon, in virtue of Lemma~\ref{LDCHEch}, we conclude that $\mathrm{ip}^{(2,X)@}_{s}(
\mathrm{CH}^{(2)}_{s}(
\mathfrak{P}^{(2)}
))$ is a   second-order echelon, thus the second-order path $\mathrm{ip}^{(2,X)@}_{s}(
\mathrm{CH}^{(2)}_{s}(
\mathfrak{Q}^{(2)}
))
\circ^{1\mathbf{Pth}_{\boldsymbol{\mathcal{A}}^{(2)}}}_{s}
\mathrm{ip}^{(2,X)@}_{s}(
\mathrm{CH}^{(2)}_{s}(
\mathfrak{P}^{(2)}
))$ is a second-order path of length strictly greater than one containing a   second-order echelon on its first step; the third equality follows from  the fact that, by Proposition~\ref{PDIpDCH}, the element $\mathrm{ip}^{(2,X)@}_{s}(
\mathrm{CH}^{(2)}_{s}(
\mathfrak{Q}^{(2)}
))$ is a second-order path in $[\mathfrak{Q}^{(2)}]^{}_{s}$ and the element $\mathrm{ip}^{(2,X)@}_{s}(
\mathrm{CH}^{(2)}_{s}(
\mathfrak{P}^{(2)}
))$ is a second-order path in $[\mathfrak{P}^{(2)}]^{}_{s}$, thus the respective values of the second-order Curry-Howard mapping are equal; finally, the last equality follows from the previous discussion on the value of the second-order Curry-Howard mapping at 
$
\mathfrak{Q}^{(2)}\circ^{1\mathbf{Pth}_{\boldsymbol{\mathcal{A}}^{(2)}}}_{s}\mathfrak{P}^{(2)}$.

The case of $\mathfrak{P}^{(2)}$ being a   second-order echelon follows.

If~(1.1.2), i.e., if $\mathfrak{Q}^{(2)}\circ^{1\mathbf{Pth}_{\boldsymbol{\mathcal{A}}^{(2)}}}_{s}\mathfrak{P}^{(2)}$ is a second-order path containing a   second-order echelon on its first step and $\mathfrak{P}^{(2)}$ is a second-order path of length strictly greater than one containing one   second-order echelon on its first step then, according to Definition~\ref{DDCH}, the value of the second-order Curry-Howard mapping at $\mathfrak{P}^{(2)}$ is given by
$$
\mathrm{CH}^{(2)}_{s}\left(
\mathfrak{P}^{(2)}
\right)
=
\mathrm{CH}^{(2)}_{s}\left(
\mathfrak{P}^{(2),1,\bb{\mathfrak{P}^{(2)}}-1}
\right)
\circ^{1\mathbf{T}_{\Sigma^{\boldsymbol{\mathcal{A}}^{(2)}}}(X)}_{s}
\mathrm{CH}^{(2)}_{s}\left(
\mathfrak{P}^{(2),0,0}
\right).
$$

Moreover, we have that $$\left(\mathfrak{Q}^{(2)}
\circ^{1\mathbf{Pth}_{\boldsymbol{\mathcal{A}}^{(2)}}}_{s}
\mathfrak{P}^{(2)}
\right)^{1,\bb{\mathfrak{Q}^{(2)}
\circ^{1\mathbf{Pth}_{\boldsymbol{\mathcal{A}}^{(2)}}}_{s}
\mathfrak{P}^{(2)}}-1}
=
\mathfrak{Q}^{(2)}
\circ^{1\mathbf{Pth}_{\boldsymbol{\mathcal{A}}^{(2)}}}_{s}
\mathfrak{P}^{(2),1,\bb{\mathfrak{P}^{(2)}}-1}
.$$

Thus, the value of the second-order Curry-Howard mapping at $\mathfrak{Q}^{(2)}
\circ^{1\mathbf{Pth}_{\boldsymbol{\mathcal{A}}^{(2)}}}_{s}
\mathfrak{P}^{(2)}$ is given by
\begin{multline*}
\mathrm{CH}^{(2)}_{s}\left(
\mathfrak{Q}^{(2)}
\circ^{1\mathbf{Pth}_{\boldsymbol{\mathcal{A}}^{(2)}}}_{s}
\mathfrak{P}^{(2)}
\right)
\\=
\mathrm{CH}^{(2)}_{s}\left(
\mathfrak{Q}^{(2)}
\circ^{1\mathbf{Pth}_{\boldsymbol{\mathcal{A}}^{(2)}}}_{s}
\mathfrak{P}^{(2),1,\bb{\mathfrak{P}^{(2)}}-1}
\right)
\circ^{1\mathbf{T}_{\Sigma^{\boldsymbol{\mathcal{A}}^{(2)}}}(X)}_{s}
\mathrm{CH}^{(2)}_{s}\left(
\mathfrak{P}^{(2),0,0}
\right).
\end{multline*}

On the other hand, the following chain of equalities holds
\begin{flushleft}
$\mathrm{CH}^{(2)}_{s}
\left(
\mathrm{ip}^{(2,X)@}_{s}
\left(
\mathrm{CH}^{(2)}_{s}\left(
\mathfrak{Q}^{(2)}
\right)
\circ^{1\mathbf{T}_{\Sigma^{\boldsymbol{\mathcal{A}}^{(2)}}}(X)}_{s}
\mathrm{CH}^{(2)}_{s}\left(
\mathfrak{P}^{(2)}
\right)
\right)
\right)
$
\allowdisplaybreaks
\begin{align*}
\quad &=
\mathrm{CH}^{(2)}_{s}
\left(
\mathrm{ip}^{(2,X)@}_{s}
\left(
\mathrm{CH}^{(2)}_{s}(
\mathfrak{Q}^{(2)}
)
\right)
\circ^{1\mathbf{Pth}_{\boldsymbol{\mathcal{A}}^{(2)}}}_{s}
\mathrm{ip}^{(2,X)@}_{s}
\left(
\mathrm{CH}^{(2)}_{s}(
\mathfrak{P}^{(2)}
)
\right)
\right)
\tag{1}
\\&=
\mathrm{CH}^{(2)}_{s}
\Bigg(
\mathrm{ip}^{(2,X)@}_{s}
\left(
\mathrm{CH}^{(2)}_{s}\left(
\mathfrak{Q}^{(2)}
\right)
\right)
\circ^{1\mathbf{Pth}_{\boldsymbol{\mathcal{A}}^{(2)}}}_{s}
\\ &\qquad
\mathrm{ip}^{(2,X)@}_{s}
\left(
\mathrm{CH}^{(2)}_{s}\left(
\mathfrak{P}^{(2),1,\bb{\mathfrak{P}^{(2)}}-1}
\right)
\circ^{1\mathbf{T}_{\Sigma^{\boldsymbol{\mathcal{A}}^{(2)}}}(X)}_{s}
\mathrm{CH}^{(2)}_{s}\left(
\mathfrak{P}^{(2),0,0}
\right)
\right)
\Bigg)
\tag{2}
\\&=
\mathrm{CH}^{(2)}_{s}
\Bigg(
\mathrm{ip}^{(2,X)@}_{s}
\left(
\mathrm{CH}^{(2)}_{s}\left(
\mathfrak{Q}^{(2)}
\right)
\right)
\circ^{1\mathbf{Pth}_{\boldsymbol{\mathcal{A}}^{(2)}}}_{s}
\\ &\qquad\qquad\qquad\qquad
\mathrm{ip}^{(2,X)@}_{s}
\left(
\mathrm{CH}^{(2)}_{s}\left(
\mathfrak{P}^{(2),1,\bb{\mathfrak{P}^{(2)}}-1}
\right)
\right)
\circ^{1\mathbf{Pth}_{\boldsymbol{\mathcal{A}}^{(2)}}}_{s}
\\ &\qquad\qquad\qquad\qquad\qquad\qquad\qquad\qquad
\mathrm{ip}^{(2,X)@}_{s}
\left(
\mathrm{CH}^{(2)}_{s}\left(
\mathfrak{P}^{(2),0,0}
\right)
\right)
\Bigg)
\tag{3}
\\&=
\mathrm{CH}^{(2)}_{s}
\left(
\mathrm{ip}^{(2,X)@}_{s}
\left(
\mathrm{CH}^{(2)}_{s}\left(
\mathfrak{Q}^{(2)}
\right)
\right)
\circ^{1\mathbf{Pth}_{\boldsymbol{\mathcal{A}}^{(2)}}}_{s}
\mathrm{ip}^{(2,X)@}_{s}
\left(
\mathrm{CH}^{(2)}_{s}\left(
\mathfrak{P}^{(2),1,\bb{\mathfrak{P}^{(2)}}-1}
\right)
\right)
\right)
\\&\qquad\qquad\qquad\qquad\quad
\circ^{1\mathbf{T}_{\Sigma^{\boldsymbol{\mathcal{A}}^{(2)}}}(X)}_{s}
\mathrm{CH}^{(2)}_{s}
\left(
\mathrm{ip}^{(2,X)@}_{s}
\left(
\mathrm{CH}^{(2)}_{s}\left(
\mathfrak{P}^{(2),0,0}
\right)
\right)
\right)
\tag{4}
\\&=
\mathrm{CH}^{(2)}_{s}
\left(
\mathrm{ip}^{(2,X)@}_{s}
\left(
\mathrm{CH}^{(2)}_{s}\left(
\mathfrak{Q}^{(2)}
\right)
\circ^{1\mathbf{T}_{\Sigma^{\boldsymbol{\mathcal{A}}^{(2)}}}(X)}_{s}
\mathrm{CH}^{(2)}_{s}\left(
\mathfrak{P}^{(2),1,\bb{\mathfrak{P}^{(2)}}-1}
\right)
\right)
\right)
\\&\qquad\qquad\qquad\qquad\qquad\qquad\qquad\qquad\qquad\quad
\circ^{1\mathbf{T}_{\Sigma^{\boldsymbol{\mathcal{A}}^{(2)}}}(X)}_{s}
\mathrm{CH}^{(2)}_{s}\left(
\mathfrak{P}^{(2),0,0}
\right)
\tag{5}
\\&=
\mathrm{CH}^{(2)}_{s}\left(
\mathfrak{Q}^{(2)}
\circ^{1\mathbf{Pth}_{\boldsymbol{\mathcal{A}}^{(2)}}}_{s}
\mathfrak{P}^{(2),1,\bb{\mathfrak{P}^{(2)}}-1}
\right)
\circ^{1\mathbf{T}_{\Sigma^{\boldsymbol{\mathcal{A}}^{(2)}}}(X)}_{s}
\mathrm{CH}^{(2)}_{s}\left(
\mathfrak{P}^{(2),0,0}
\right)
\tag{6}
\\&=
\mathrm{CH}^{(2)}_{s}\left(
\mathfrak{Q}^{(2)}
\circ^{1\mathbf{Pth}_{\boldsymbol{\mathcal{A}}^{(2)}}}_{s}
\mathfrak{P}^{(2)}
\right).
\tag{7}
\end{align*}
\end{flushleft}

The first equality is a direct consequence of the equality proved in the first item of this lemma; the second equality unravels the definition of the second-order Curry-Howard mapping at $\mathfrak{P}^{(2)}$; the third equality follows from the fact that, by Definition~\ref{DDIp}, $\mathrm{ip}^{(2,X)@}$ is a $\Sigma^{\boldsymbol{\mathcal{A}}^{(2)}}$-homomorphism and from the fact that, by Proposition~\ref{PDIpDCH},  $\mathrm{ip}^{(2,X)@}_{s}(
\mathrm{CH}^{(2)}_{s}(
\mathfrak{P}^{(2),1,\bb{\mathfrak{P}^{(2)}}-1}
))$ and $\mathrm{ip}^{(2,X)@}_{s}(
\mathrm{CH}^{(2)}_{s}(
\mathfrak{P}^{(2),0,0}
))$ are second-order paths in $\mathrm{Pth}_{\boldsymbol{\mathcal{A}}^{(2)},s}$. We have omitted parenthesis because, by Proposition~\ref{PDPthComp}, the $1$-composition of second-order paths is associative; the fourth equality follows from the fact that $\mathfrak{P}^{(2),0,0}$ is a   second-order echelon and, by  Proposition~\ref{PDIpDCH}, $\mathrm{ip}^{(2,X)@}_{s}(\mathrm{CH}^{(2)}_{s}(\mathfrak{P}^{(2),0,0}))$ is a second-order path in $[\mathfrak{P}^{(2),0,0}]^{}_{s}$ then, by Lemma~\ref{LDCHEch}, we have that $\mathrm{ip}^{(2,X)@}_{s}(
\mathrm{CH}^{(2)}_{s}(
\mathfrak{P}^{(2),0,0}
))$ is a   second-order echelon. Also, by Proposition~\ref{PDIpDCH}, $\mathrm{ip}^{(2,X)@}_{s}(\mathrm{CH}^{(2)}_{s}(\mathfrak{P}^{(2),1,\bb{\mathfrak{P}^{(2)}}-1}))$ is a second-order path in $[\mathfrak{P}^{(2),1,\bb{\mathfrak{P}^{(2)}}-1}]^{}_{s}$ then, by Lemma~\ref{LDCH}, we have that the second-order paths $\mathrm{ip}^{(2,X)@}_{s}(\mathrm{CH}^{(2)}_{s}(\mathfrak{P}^{(2),1,\bb{\mathfrak{P}^{(2)}}-1}))$ and $\mathfrak{P}^{(2),1,\bb{\mathfrak{P}^{(2)}}-1}$ have the same length. So, considering the foregoing, we have that the second-order path under consideration is a second-order path of length strictly greater than one containing a   second-order echelon on its first step. The equality simply unravels the second-order Curry-Howard mapping at a second-order path of such nature; the fifth equality follows, on the one hand, from the fact that, by Definition~\ref{DDIp}, $\mathrm{ip}^{(2,X)@}$ is a $\Sigma^{\boldsymbol{\mathcal{A}}^{(2)}}$-homomorphism and, on the other hand, from the fact that, by Proposition~\ref{PDIpDCH}, $\mathrm{ip}^{(2,X)@}_{s}(\mathrm{CH}^{(2)}_{s}(\mathfrak{P}^{(2),0,0}))$ is a second-order path in $[\mathfrak{P}^{(2),0,0}]^{}_{s}$; the sixth equality follows from induction. Let us note that the pair $(\mathfrak{Q}^{(2)}\circ^{1\mathbf{Pth}_{\boldsymbol{\mathcal{A}}^{(2)}}}_{s}\mathfrak{P}^{(2),1,\bb{\mathfrak{P}^{(2)}}-1},s)$ $\prec_{\mathbf{Pth}_{\boldsymbol{\mathcal{A}}^{(2)}}}$-precedes $(\mathfrak{Q}^{(2)}\circ^{1\mathbf{Pth}_{\boldsymbol{\mathcal{A}}^{(2)}}}_{s}\mathfrak{P}^{(2)},s)$, thus we have that
\begin{multline*}
\mathrm{CH}^{(2)}_{s}
\left(
\mathrm{ip}^{(2,X)@}_{s}
\left(
\mathrm{CH}^{(2)}_{s}\left(
\mathfrak{Q}^{(2)}
\right)
\circ^{1\mathbf{T}_{\Sigma^{\boldsymbol{\mathcal{A}}^{(2)}}}(X)}_{s}
\mathrm{CH}^{(2)}_{s}\left(
\mathfrak{P}^{(2),1,\bb{\mathfrak{P}^{(2)}}-1}
\right)
\right)
\right)\\ =
\mathrm{CH}^{(2)}_{s}\left(
\mathfrak{Q}^{(2)}
\circ^{1\mathbf{Pth}_{\boldsymbol{\mathcal{A}}^{(2)}}}_{s}
\mathfrak{P}^{(2),1,\bb{\mathfrak{P}^{(2)}}-1}
\right);
\end{multline*}
finally, the last equality recovers the value of the second-order Curry-Howard mapping at $\mathfrak{Q}^{(2)}\circ^{1\mathbf{Pth}_{\boldsymbol{\mathcal{A}}^{(2)}}}_{s}\mathfrak{P}^{(2)}$.

The case of $\mathfrak{P}^{(2)}$ being a second-order path of length strictly greater than one containing a   second-order echelon on its first step follows.

This conclude the case $i=0$.

Now, if~(1.2), i.e., if $\mathfrak{Q}^{(2)}\circ^{1\mathbf{Pth}_{\boldsymbol{\mathcal{A}}^{(2)}}}_{s}\mathfrak{P}^{(2)}$ is a second-order path of length strictly greater than one containing its first   echelon at position $i\in\bb{\mathfrak{Q}^{(2)}\circ^{1\mathbf{Pth}_{\boldsymbol{\mathcal{A}}^{(2)}}}_{s}\mathfrak{P}^{(2)}}$ with $i>0$, then since $\bb{\mathfrak{Q}^{(2)}\circ^{1\mathbf{Pth}_{\boldsymbol{\mathcal{A}}^{(2)}}}_{s}\mathfrak{P}^{(2)}}=\bb{\mathfrak{Q}^{(2)}}+\bb{\mathfrak{P}^{(2)}}$, it could be the case that either (1.2.1) $i\in\bb{\mathfrak{P}^{(2)}}$ or (1.2.2) $i\in[
\bb{\mathfrak{P}^{(2)}},\bb{\mathfrak{P}^{(2)}}+\bb{\mathfrak{Q}^{(2)}}-1
]$.

If~(1.2.1), i.e., if $\mathfrak{Q}^{(2)}\circ^{1\mathbf{Pth}_{\boldsymbol{\mathcal{A}}^{(2)}}}_{s}\mathfrak{P}^{(2)}$ is a second-order path of length strictly greater than one containing its first   echelon at position $i\in\bb{\mathfrak{P}^{(2)}}$ with $i>0$, then $\mathfrak{P}^{(2)}$ is a second-order path of length strictly greater than one containing a   second-order echelon on a step different from the initial one. Then, according to Definition~\ref{DDCH}, the value of the second-order Curry-Howard mapping at $\mathfrak{P}^{(2)}$ is given by
$$
\mathrm{CH}^{(2)}_{s}\left(
\mathfrak{P}^{(2)}
\right)
=
\mathrm{CH}^{(2)}_{s}\left(
\mathfrak{P}^{(2),i,\bb{\mathfrak{P}^{(2)}}-1}
\right)
\circ^{1\mathbf{T}_{\Sigma^{\boldsymbol{\mathcal{A}}^{(2)}}}(X)}_{s}
\mathrm{CH}^{(2)}_{s}\left(
\mathfrak{P}^{(2),0,i-1}
\right).
$$

Moreover, we have that $$\left(\mathfrak{Q}^{(2)}
\circ^{1\mathbf{Pth}_{\boldsymbol{\mathcal{A}}^{(2)}}}_{s}
\mathfrak{P}^{(2)}
\right)^{i,\bb{\mathfrak{Q}^{(2)}
\circ^{1\mathbf{Pth}_{\boldsymbol{\mathcal{A}}^{(2)}}}_{s}
\mathfrak{P}^{(2)}}-1}
=
\mathfrak{Q}^{(2)}
\circ^{1\mathbf{Pth}_{\boldsymbol{\mathcal{A}}^{(2)}}}_{s}
\mathfrak{P}^{(2),i,\bb{\mathfrak{P}^{(2)}}-1}
.$$

Thus, the value of the second-order Curry-Howard mapping at $\mathfrak{Q}^{(2)}
\circ^{1\mathbf{Pth}_{\boldsymbol{\mathcal{A}}^{(2)}}}_{s}
\mathfrak{P}^{(2)}$ is given by
\begin{multline*}
\mathrm{CH}^{(2)}_{s}\left(
\mathfrak{Q}^{(2)}
\circ^{1\mathbf{Pth}_{\boldsymbol{\mathcal{A}}^{(2)}}}_{s}
\mathfrak{P}^{(2)}
\right)
\\=
\mathrm{CH}^{(2)}_{s}\left(
\mathfrak{Q}^{(2)}
\circ^{1\mathbf{Pth}_{\boldsymbol{\mathcal{A}}^{(2)}}}_{s}
\mathfrak{P}^{(2),i,\bb{\mathfrak{P}^{(2)}}-1}
\right)
\circ^{1\mathbf{T}_{\Sigma^{\boldsymbol{\mathcal{A}}^{(2)}}}(X)}_{s}
\mathrm{CH}^{(2)}_{s}\left(
\mathfrak{P}^{(2),0,i-1}
\right).
\end{multline*}

On the other hand, the following chain of equalities holds
\begin{flushleft}
$\mathrm{CH}^{(2)}_{s}
\left(
\mathrm{ip}^{(2,X)@}_{s}
\left(
\mathrm{CH}^{(2)}_{s}\left(
\mathfrak{Q}^{(2)}
\right)
\circ^{1\mathbf{T}_{\Sigma^{\boldsymbol{\mathcal{A}}^{(2)}}}(X)}_{s}
\mathrm{CH}^{(2)}_{s}\left(
\mathfrak{P}^{(2)}
\right)
\right)
\right)
$
\allowdisplaybreaks
\begin{align*}
\quad &=
\mathrm{CH}^{(2)}_{s}
\left(
\mathrm{ip}^{(2,X)@}_{s}
\left(
\mathrm{CH}^{(2)}_{s}\left(
\mathfrak{Q}^{(2)}
\right)
\right)
\circ^{1\mathbf{Pth}_{\boldsymbol{\mathcal{A}}^{(2)}}}_{s}
\mathrm{ip}^{(2,X)@}_{s}
\left(
\mathrm{CH}^{(2)}_{s}\left(
\mathfrak{P}^{(2)}
\right)
\right)
\right)
\tag{1}
\\&=
\mathrm{CH}^{(2)}_{s}
\Bigg(
\mathrm{ip}^{(2,X)@}_{s}
\left(
\mathrm{CH}^{(2)}_{s}\left(
\mathfrak{Q}^{(2)}
\right)
\right)
\circ^{1\mathbf{Pth}_{\boldsymbol{\mathcal{A}}^{(2)}}}_{s}
\\ &\qquad\quad
\mathrm{ip}^{(2,X)@}_{s}
\left(
\mathrm{CH}^{(2)}_{s}\left(
\mathfrak{P}^{(2),i,\bb{\mathfrak{P}^{(2)}}-1}
\right)
\circ^{1\mathbf{T}_{\Sigma^{\boldsymbol{\mathcal{A}}^{(2)}}}(X)}_{s}
\mathrm{CH}^{(2)}_{s}\left(
\mathfrak{P}^{(2),0,i-1}
\right)
\right)
\Bigg)
\tag{2}
\\&=
\mathrm{CH}^{(2)}_{s}
\Bigg(
\mathrm{ip}^{(2,X)@}_{s}
\left(
\mathrm{CH}^{(2)}_{s}\left(
\mathfrak{Q}^{(2)}
\right)
\right)
\circ^{1\mathbf{Pth}_{\boldsymbol{\mathcal{A}}^{(2)}}}_{s}
\\&\qquad\qquad\qquad\qquad
\mathrm{ip}^{(2,X)@}_{s}
\left(
\mathrm{CH}^{(2)}_{s}\left(
\mathfrak{P}^{(2),i,\bb{\mathfrak{P}^{(2)}}-1}
\right)
\right)
\circ^{1\mathbf{Pth}_{\boldsymbol{\mathcal{A}}^{(2)}}}_{s}
\\&\qquad\qquad\qquad\qquad\qquad\qquad\qquad\qquad
\mathrm{ip}^{(2,X)@}_{s}
\left(
\mathrm{CH}^{(2)}_{s}\left(
\mathfrak{P}^{(2),0,i-1}
\right)
\right)
\Bigg)
\tag{3}
\\&=
\mathrm{CH}^{(2)}_{s}
\left(
\mathrm{ip}^{(2,X)@}_{s}
\left(
\mathrm{CH}^{(2)}_{s}\left(
\mathfrak{Q}^{(2)}
\right)
\right)
\circ^{1\mathbf{Pth}_{\boldsymbol{\mathcal{A}}^{(2)}}}_{s}
\mathrm{ip}^{(2,X)@}_{s}
\left(
\mathrm{CH}^{(2)}_{s}\left(
\mathfrak{P}^{(2),i,\bb{\mathfrak{P}^{(2)}}-1}
\right)
\right)
\right)
\\&\qquad\qquad
\circ^{1\mathbf{T}_{\Sigma^{\boldsymbol{\mathcal{A}}^{(2)}}}(X)}_{s}
\mathrm{CH}^{(2)}_{s}
\left(
\mathrm{ip}^{(2,X)@}_{s}
\left(
\mathrm{CH}^{(2)}_{s}\left(
\mathfrak{P}^{(2),0,i-1}
\right)
\right)
\right)
\tag{4}
\\&=
\mathrm{CH}^{(2)}_{s}
\left(
\mathrm{ip}^{(2,X)@}_{s}
\left(
\mathrm{CH}^{(2)}_{s}\left(
\mathfrak{Q}^{(2)}
\right)
\circ^{1\mathbf{T}_{\Sigma^{\boldsymbol{\mathcal{A}}^{(2)}}}(X)}_{s}
\mathrm{CH}^{(2)}_{s}\left(
\mathfrak{P}^{(2),i,\bb{\mathfrak{P}^{(2)}}-1}
\right)
\right)
\right)
\\&\qquad\qquad\qquad\qquad\qquad\qquad\qquad\qquad\qquad
\circ^{1\mathbf{T}_{\Sigma^{\boldsymbol{\mathcal{A}}^{(2)}}}(X)}_{s}
\mathrm{CH}^{(2)}_{s}\left(
\mathfrak{P}^{(2),0,i-1}
\right)
\tag{5}
\\&=
\mathrm{CH}^{(2)}_{s}\left(
\mathfrak{Q}^{(2)}
\circ^{1\mathbf{Pth}_{\boldsymbol{\mathcal{A}}^{(2)}}}_{s}
\mathfrak{P}^{(2),i,\bb{\mathfrak{P}^{(2)}}-1}
\right)
\circ^{1\mathbf{T}_{\Sigma^{\boldsymbol{\mathcal{A}}^{(2)}}}(X)}_{s}
\mathrm{CH}^{(2)}_{s}\left(
\mathfrak{P}^{(2),0,i-1}
\right)
\tag{6}
\\&=
\mathrm{CH}^{(2)}_{s}\left(
\mathfrak{Q}^{(2)}
\circ^{1\mathbf{Pth}_{\boldsymbol{\mathcal{A}}^{(2)}}}_{s}
\mathfrak{P}^{(2)}
\right).
\tag{7}
\end{align*}
\end{flushleft}

The first equality is a direct consequence of the equality proved in the first item of this lemma; the second equality unravels the definition of the second-order Curry-Howard mapping at $\mathfrak{P}^{(2)}$; the third equality follows from the fact that, by Definition~\ref{DDIp}, $\mathrm{ip}^{(2,X)@}$ is a $\Sigma^{\boldsymbol{\mathcal{A}}^{(2)}}$-homomorphism and from the fact that, by Proposition~\ref{PDIpDCH}, the elements $\mathrm{ip}^{(2,X)@}_{s}(
\mathrm{CH}^{(2)}_{s}(
\mathfrak{P}^{(2),i,\bb{\mathfrak{P}^{(2)}}-1}
))$ and $\mathrm{ip}^{(2,X)@}_{s}(
\mathrm{CH}^{(2)}_{s}(
\mathfrak{P}^{(2),0,i-1}
))$ are second-order paths in $\mathrm{Pth}_{\boldsymbol{\mathcal{A}}^{(2)},s}$. We have omitted parentheses because, by Proposition~\ref{PDPthComp}, the $1$-composition of second-order paths is associative;  the fourth equality follows from the fact that $\mathfrak{P}^{(2),0,i-1}$ is an echelonless second-order path and, by Proposition~\ref{PDIpDCH}, $\mathrm{ip}^{(2,X)@}_{s}(\mathrm{CH}^{(2)}_{s}(\mathfrak{P}^{(2),0,i-1}))$ is a second-order path in $[\mathfrak{P}^{(2),0,i-1}]^{}_{s}$ then, by Lemma~\ref{LDCHNEch}, we see that $\mathrm{ip}^{(2,X)@}_{s}(
\mathrm{CH}^{(2)}_{s}(
\mathfrak{P}^{(2),0,i-1}
))$ is an echelonless second-order path. By Proposition~\ref{PDIpDCH}, we also have that  $\mathrm{ip}^{(2,X)@}_{s}(\mathrm{CH}^{(2)}_{s}(\mathfrak{P}^{(2),i,\bb{\mathfrak{P}^{(2)}}-1}))$ is a second-order path in $[\mathfrak{P}^{(2),i,\bb{\mathfrak{P}^{(2)}}-1}]^{}_{s}$ then, by Lemma~\ref{LDCHEchInt}, we have that $\mathrm{ip}^{(2,X)@}_{s}(\mathrm{CH}^{(2)}_{s}(\mathfrak{P}^{(2),i,\bb{\mathfrak{P}^{(2)}}-1}))$ is a second-order path containing a   second-order echelon on its first step.  All things considered, we have that the second-order path under consideration is a second-order path of length strictly greater than one containing a   second-order echelon at position $i$. The equality simply unravels the second-order Curry-Howard mapping at a second-order path of such nature; the fifth equality follows, on the one hand, from the fact that, by Definition~\ref{DDIp}, $\mathrm{ip}^{(2,X)@}$ is a $\Sigma^{\boldsymbol{\mathcal{A}}^{(2)}}$-homomorphism and, on the other hand, from the fact that, by Proposition~\ref{PDIpDCH}, $\mathrm{ip}^{(2,X)@}_{s}(\mathrm{CH}^{(2)}_{s}(\mathfrak{P}^{(2),0,i-1}))$ is a second-order path in $[\mathfrak{P}^{(2),0,i-1}]^{}_{s}$; the sixth equality follows from induction. Let us note that the pair $(\mathfrak{Q}^{(2)}\circ^{1\mathbf{Pth}_{\boldsymbol{\mathcal{A}}^{(2)}}}_{s}\mathfrak{P}^{(2),i,\bb{\mathfrak{P}^{(2)}}-1},s)$ $\prec_{\mathbf{Pth}_{\boldsymbol{\mathcal{A}}^{(2)}}}$-precedes $(\mathfrak{Q}^{(2)}\circ^{1\mathbf{Pth}_{\boldsymbol{\mathcal{A}}^{(2)}}}_{s}\mathfrak{P}^{(2)},s)$, thus we have that
\begin{multline*}
\mathrm{CH}^{(2)}_{s}
\left(
\mathrm{ip}^{(2,X)@}_{s}
\left(
\mathrm{CH}^{(2)}_{s}\left(
\mathfrak{Q}^{(2)}
\right)
\circ^{1\mathbf{T}_{\Sigma^{\boldsymbol{\mathcal{A}}^{(2)}}}(X)}_{s}
\mathrm{CH}^{(2)}_{s}\left(
\mathfrak{P}^{(2),i,\bb{\mathfrak{P}^{(2)}}-1}
\right)
\right)
\right)\\ =
\mathrm{CH}^{(2)}_{s}\left(
\mathfrak{Q}^{(2)}
\circ^{1\mathbf{Pth}_{\boldsymbol{\mathcal{A}}^{(2)}}}_{s}
\mathfrak{P}^{(2),i,\bb{\mathfrak{P}^{(2)}}-1}
\right);
\end{multline*}
finally, the last equality recovers the value of the second-order Curry-Howard mapping at $\mathfrak{Q}^{(2)}\circ^{1\mathbf{Pth}_{\boldsymbol{\mathcal{A}}^{(2)}}}_{s}\mathfrak{P}^{(2)}$.

The case $i\in\bb{\mathfrak{P}^{(2)}}$ follows.

If~(1.2.2), i.e., if $\mathfrak{Q}^{(2)}\circ^{1\mathbf{Pth}_{\boldsymbol{\mathcal{A}}^{(2)}}}_{s}\mathfrak{P}^{(2)}$ is a second-order path of length strictly greater than one containing its first   echelon at position $i\in[\bb{\mathfrak{P}^{(2)}}, \bb{\mathfrak{P}^{(2)}}+\bb{\mathfrak{Q}^{(2)}}-1]$, then $\mathfrak{Q}^{(2)}$ is a non-$(2,[1])$-identity second-order path containing a   second-order echelon, whilst $\mathfrak{P}^{(2)}$ is an echelonless second-order path.

We will distinguish three cases according to whether (1.2.2.1) $\mathfrak{Q}^{(2)}$ is a   second-order echelon; (1.2.2.2) $\mathfrak{Q}^{(2)}$ is a second-order path of length strictly greater than one containing a   second-order echelon on its first step or (1.2.2.3) $\mathfrak{Q}^{(2)}$ is a second-order path of length strictly greater than one containing a   second-order echelon on a step different from zero. These cases can be proven following a similar argument to those three cases presented above. We leave the details for the interested reader.

If~(2), i.e., if $\mathfrak{Q}^{(2)}\circ^{1\mathbf{Pth}_{\boldsymbol{\mathcal{A}}^{(2)}}}_{s}
\mathfrak{P}^{(2)}$ is an echelonless second-order path, it could be the case that (2.1) $\mathfrak{Q}^{(2)}\circ^{1\mathbf{Pth}_{\boldsymbol{\mathcal{A}}^{(2)}}}_{s}
\mathfrak{P}^{(2)}$ is an echelonless second-order path that is not head-constant, or (2.2) $\mathfrak{Q}^{(2)}\circ^{1\mathbf{Pth}_{\boldsymbol{\mathcal{A}}^{(2)}}}_{s}
\mathfrak{P}^{(2)}$ is a head-constant echelonless second-order path that is not coherent, or (2.3) $\mathfrak{Q}^{(2)}\circ^{1\mathbf{Pth}_{\boldsymbol{\mathcal{A}}^{(2)}}}_{s}
\mathfrak{P}^{(2)}$ is a coherent head-constant echelonless second-order path.

If (2.1), i.e., if $\mathfrak{Q}^{(2)}\circ^{1\mathbf{Pth}_{\boldsymbol{\mathcal{A}}^{(2)}}}_{s}
\mathfrak{P}^{(2)}$ is an echelonless second-order path that is not head-constant we let $i\in\bb{\mathfrak{Q}^{(2)}\circ^{1\mathbf{Pth}_{\boldsymbol{\mathcal{A}}^{(2)}}}_{s}
\mathfrak{P}^{(2)}}$ be the greatest index for which $(\mathfrak{Q}^{(2)}\circ^{1\mathbf{Pth}_{\boldsymbol{\mathcal{A}}^{(2)}}}_{s}
\mathfrak{P}^{(2)})^{0,i}$ is a head-constant echelonless second-order path. Since $\bb{\mathfrak{Q}^{(2)}\circ^{1\mathbf{Pth}_{\boldsymbol{\mathcal{A}}^{(2)}}}_{s}
\mathfrak{P}^{(2)}}=\bb{\mathfrak{Q}^{(2)}}+\bb{\mathfrak{P}^{(2)}}$, then either (2.1.1) $i\in\bb{\mathfrak{P}^{(2)}}-1$, (2.1.2) $i=\bb{\mathfrak{P}^{(2)}}-1$, or (2.1.3) $i\in [
\bb{\mathfrak{P}^{(2)}}, \bb{\mathfrak{Q}^{(2)}
\circ^{1\mathbf{Pth}_{\boldsymbol{\mathcal{A}}^{(2)}}}_{s}
\mathfrak{P}^{(2)}}-1
]$.

If~(2.1.1), i.e.,  if $\mathfrak{Q}^{(2)}\circ^{1\mathbf{Pth}_{\boldsymbol{\mathcal{A}}^{(2)}}}_{s}
\mathfrak{P}^{(2)}$ is an echelonless second-order path that is not head-constant and $i\in\bb{\mathfrak{P}^{(2)}}$ is the greatest index for which $(\mathfrak{Q}^{(2)}\circ^{1\mathbf{Pth}_{\boldsymbol{\mathcal{A}}^{(2)}}}_{s}
\mathfrak{P}^{(2)})^{0,i}$ is a head-constant echelonless second-order path then, regarding the second-order paths $\mathfrak{Q}^{(2)}$ and $\mathfrak{P}^{(2)}$, we have that
\begin{itemize}
\item[(i)] $\mathfrak{P}^{(2),0,i}$ is a coherent  head-constant echelonless second-order path; Moreover, $i\in\bb{\mathfrak{P}^{(2)}}-1$ is the greatest index for which $\mathfrak{P}^{(2),0,i}$ is head-constant;
\item[(ii)] $\mathfrak{P}^{(2),i+1,\bb{\mathfrak{P}^{(2)}}-1}$ is an echelonless second-order path;
\item[(iii)] $\mathfrak{Q}^{(2)}$ is an echelonless second-order path.
\end{itemize}

Following (i) and (ii) and taking into account Definition~\ref{DDCH}, we have that the value of the second-order Curry-Howard mapping at $\mathfrak{P}^{(2)}$ is given by
$$
\mathrm{CH}^{(2)}_{s}\left(
\mathfrak{P}^{(2)}
\right)
=
\mathrm{CH}^{(2)}_{s}\left(
\mathfrak{P}^{(2),i+1,\bb{\mathfrak{P}^{(2)}}-1}
\right)
\circ^{1\mathbf{T}_{\Sigma^{\boldsymbol{\mathcal{A}}^{(2)}}}}_{s}
\mathrm{CH}^{(2)}_{s}\left(
\mathfrak{P}^{(2),0,i}
\right).
$$

Moreover, we have that 
$$
\left(\mathfrak{Q}^{(2)}\circ^{1\mathbf{Pth}_{\boldsymbol{\mathcal{A}}^{(2)}}}_{s}
\mathfrak{P}^{(2)}\right)^{i+1,\bb{\mathfrak{Q}^{(2)}\circ^{1\mathbf{Pth}_{\boldsymbol{\mathcal{A}}^{(2)}}}_{s}
\mathfrak{P}^{(2)}}-1}
=
\mathfrak{Q}^{(2)}
\circ^{1\mathbf{Pth}_{\boldsymbol{\mathcal{A}}^{(2)}}}_{s}
\mathfrak{P}^{(2),i+1,\bb{\mathfrak{P}^{(2)}}-1}.
$$

Thus, the value of the second-order Curry-Howard mapping at $\mathfrak{Q}^{(2)}\circ^{1\mathbf{Pth}_{\boldsymbol{\mathcal{A}}^{(2)}}}_{s}
\mathfrak{P}^{(2)}$ is given by
\begin{multline*}
\mathrm{CH}^{(2)}_{s}\left(
\mathfrak{Q}^{(2)}\circ^{1\mathbf{Pth}_{\boldsymbol{\mathcal{A}}^{(2)}}}_{s}
\mathfrak{P}^{(2)}
\right)
\\=
\mathrm{CH}^{(2)}_{s}\left(
\mathfrak{Q}^{(2)}
\circ^{1\mathbf{Pth}_{\boldsymbol{\mathcal{A}}^{(2)}}}_{s}
\mathfrak{P}^{(2),i+1,\bb{\mathfrak{P}^{(2)}}-1}
\right)
\circ^{1\mathbf{T}_{\Sigma^{\boldsymbol{\mathcal{A}}^{(2)}}}(X)}_{s}
\mathrm{CH}^{(2)}_{s}\left(
\mathfrak{P}^{(2),0,i}
\right).
\end{multline*}

On the other hand, the following chain of equalities holds
\begin{flushleft}
$\mathrm{CH}^{(2)}_{s}
\left(
\mathrm{ip}^{(2,X)@}_{s}
\left(
\mathrm{CH}^{(2)}_{s}\left(
\mathfrak{Q}^{(2)}
\right)
\circ^{1\mathbf{T}_{\Sigma^{\boldsymbol{\mathcal{A}}^{(2)}}}(X)}_{s}
\mathrm{CH}^{(2)}_{s}\left(
\mathfrak{P}^{(2)}
\right)
\right)
\right)
$
\allowdisplaybreaks
\begin{align*}
\quad&=
\mathrm{CH}^{(2)}_{s}
\left(
\mathrm{ip}^{(2,X)@}_{s}
\left(
\mathrm{CH}^{(2)}_{s}\left(
\mathfrak{Q}^{(2)}
\right)
\right)
\circ^{1\mathbf{Pth}_{\boldsymbol{\mathcal{A}}^{(2)}}}_{s}
\mathrm{ip}^{(2,X)@}_{s}
\left(
\mathrm{CH}^{(2)}_{s}\left(
\mathfrak{P}^{(2)}
\right)
\right)
\right)
\tag{1}
\\&=
\mathrm{CH}^{(2)}_{s}
\Bigg(
\mathrm{ip}^{(2,X)@}_{s}
\left(
\mathrm{CH}^{(2)}_{s}\left(
\mathfrak{Q}^{(2)}
\right)
\right)
\circ^{1\mathbf{Pth}_{\boldsymbol{\mathcal{A}}^{(2)}}}_{s}
\\&\qquad
\mathrm{ip}^{(2,X)@}_{s}
\left(
\mathrm{CH}^{(2)}_{s}\left(
\mathfrak{P}^{(2),i+1,\bb{\mathfrak{P}^{(2)}}-1}
\right)
\circ^{1\mathbf{T}_{\Sigma^{\boldsymbol{\mathcal{A}}^{(2)}}}(X)}_{s}
\mathrm{CH}^{(2)}_{s}\left(
\mathfrak{P}^{(2),0,i}
\right)
\right)
\Bigg)
\tag{2}
\\&=
\mathrm{CH}^{(2)}_{s}
\Bigg(
\mathrm{ip}^{(2,X)@}_{s}
\left(
\mathrm{CH}^{(2)}_{s}\left(
\mathfrak{Q}^{(2)}
\right)
\right)
\circ^{1\mathbf{Pth}_{\boldsymbol{\mathcal{A}}^{(2)}}}_{s}
\\&\qquad\qquad\qquad\qquad
\mathrm{ip}^{(2,X)@}_{s}
\left(
\mathrm{CH}^{(2)}_{s}\left(
\mathfrak{P}^{(2),i+1,\bb{\mathfrak{P}^{(2)}}-1}
\right)
\right)
\circ^{1\mathbf{Pth}_{\boldsymbol{\mathcal{A}}^{(2)}}}_{s}
\\&\qquad\qquad\qquad\qquad\qquad\qquad\qquad\qquad
\mathrm{ip}^{(2,X)@}_{s}
\left(
\mathrm{CH}^{(2)}_{s}\left(
\mathfrak{P}^{(2),0,i}
\right)
\right)
\Bigg).
\tag{3}
\end{align*}
\end{flushleft}

The first equality is a direct consequence of the equation proven in the first item of this lemma; the second equality unravels the value of the second-order Curry-Howard mapping at $\mathfrak{P}^{(2)}$; the third equality follows from the fact that, by Definition~\ref{DDIp}, $\mathrm{ip}^{(2,X)@}$ is a $\Sigma^{\boldsymbol{\mathcal{A}}^{(2)}}$-homomorphism and from the fact that, by Proposition~\ref{PDIpDCH}, the elements $\mathrm{ip}^{(2,X)@}_{s}(\mathrm{CH}^{(2)}_{s}(\mathfrak{P}^{(2),i,\bb{\mathfrak{P}^{(2)}}-1}))$ and $\mathrm{ip}^{(2,X)@}_{s}(\mathrm{CH}^{(2)}_{s}(\mathfrak{P}^{(2),0,i}))$ are second-order paths in $\mathrm{Pth}_{\boldsymbol{\mathcal{A}}^{(2)},s}$. We have omitted parenthesis because, by Proposition~\ref{PDPthComp}, the $1$-composition of second-order paths is associative.

We will center our efforts in studying the nature of the second-order path $\mathfrak{R}^{(2)}$, given by
\begin{multline*}
\mathfrak{R}^{(2)}=
\mathrm{ip}^{(2,X)@}_{s}
\left(
\mathrm{CH}^{(2)}_{s}\left(
\mathfrak{Q}^{(2)}
\right)
\right)
\circ^{1\mathbf{Pth}_{\boldsymbol{\mathcal{A}}^{(2)}}}_{s}
\\
\mathrm{ip}^{(2,X)@}_{s}
\left(
\mathrm{CH}^{(2)}_{s}\left(
\mathfrak{P}^{(2),i+1,\bb{\mathfrak{P}^{(2)}}-1}
\right)
\right)
\circ^{1\mathbf{Pth}_{\boldsymbol{\mathcal{A}}^{(2)}}}_{s}
\\
\mathrm{ip}^{(2,X)@}_{s}
\left(
\mathrm{CH}^{(2)}_{s}\left(
\mathfrak{P}^{(2),0,i}
\right)
\right)
\end{multline*}
Note that
\begin{itemize}
\item[(i)]  $\mathfrak{P}^{(2),0,i}$ is a head-constant echelonless second-order path. 
\end{itemize}

Therefore, for a unique word $\mathbf{s}\in S^{\star}-\{\lambda\}$ and a unique operation symbol $\tau\in\Sigma^{\boldsymbol{\mathcal{A}}}_{\mathbf{s},s}$, the family of first-order translations of occurring in $\mathfrak{P}^{(2),0,i}$ is a family of first-order translations of type $\tau$.

By Lemma~\ref{LDCHNEchHd}, we have that $\mathrm{CH}^{(2)}_{s}(\mathfrak{P}^{(2),0,i})\in\mathcal{T}(\tau,\mathrm{T}_{\Sigma^{\boldsymbol{\mathcal{A}}^{(2)}}}(X))^{\star}$, which is a subset of $[\mathrm{T}_{\Sigma^{\boldsymbol{\mathcal{A}}^{(2)}}}(X)]^{\mathsf{HdC}}_{s}$. Therefore, since, by Proposition~\ref{PDIpDCH}, $\mathrm{ip}^{(2,X)@}_{s}(\mathrm{CH}^{(2)}_{s}(\mathfrak{P}^{(2),0,i}))$ is a second-order path in $[\mathfrak{P}^{(2),0,i}]^{}_{s}$, then, by Lemma~\ref{LDCHNEchHd}, we have that
\begin{itemize}
\item[(i)] $\mathrm{ip}^{(2,X)@}_{s}(\mathrm{CH}^{(2)}_{s}(\mathfrak{P}^{(2),0,i}))$ is a head-constant echelonless second-order path  associated to the operation symbol $\tau$.
\end{itemize}

On the other hand 
\begin{itemize}
\item[(ii)] $\mathfrak{P}^{(2),i+1,\bb{\mathfrak{P}^{(2)}}-1}$ is an echelonless second-order path.
\end{itemize}

Therefore, for a unique word $\mathbf{s}'\in S^{\star}-\{\lambda\}$ and a unique operation symbol $\tau'\in\Sigma^{\boldsymbol{\mathcal{A}}}_{\mathbf{s}',s}$ the first  first-order translation of occurring in $\mathfrak{P}^{(2),i+1,\bb{\mathfrak{P}^{(2)}}-1}$ is a   first-order translation of type $\tau'$.

Let us recall that $\tau\neq\tau'$ since $\mathfrak{P}^{(2)}$ is an echelonless second-order path that is not head-constant for which $i\in\bb{\mathfrak{P}^{(2)}}-1$ is the greatest index for which $\mathfrak{P}^{(2),0,i}$ is a head-constant echelonless second-order path. 

By Lemma~\ref{LDCHNEch}, we have that 
$$\mathrm{CH}^{(2)}_{s}\left(\mathfrak{P}^{(2),i+1,\bb{\mathfrak{P}^{(2)}}-1}\right)\in\mathrm{T}_{\Sigma^{\boldsymbol{\mathcal{A}}^{(2)}}}(X)_{s}
-\left(
\eta^{(2,1)\sharp}[
\mathrm{PT}_{\boldsymbol{\mathcal{A}}}
]_{s}
\cup
\eta^{(2,\mathcal{A}^{(2)})}\left[\mathcal{A}^{(2)}\right]^{\mathrm{pct}}_{s}
\right).
$$
Therefore, by Proposition~\ref{PDIpDCH}, $\mathrm{ip}^{(2,X)@}_{s}(\mathrm{CH}^{(2)}_{s}(\mathfrak{P}^{(2),i+1,\bb{\mathfrak{P}^{(2)}}-1}))$ is a second-order path in $[\mathfrak{P}^{(2),i+1,\bb{\mathfrak{P}^{(2)}}-1}]^{}_{s}$, then, by Lemma~\ref{LDCHNEch}, we have that 
\begin{itemize}
\item[(ii)] $\mathrm{ip}^{(2,X)@}_{s}(\mathrm{CH}^{(2)}_{s}(\mathfrak{P}^{(2),i+1,\bb{\mathfrak{P}^{(2)}}-1}))$ is an echelonless second-order path whose first first-order translations is associated to the operation symbol $\tau'$.
\end{itemize}

Finally, 
\begin{itemize}
\item[(iii)] $\mathfrak{Q}^{(2)}$ is an echelonless second-order path.
\end{itemize}

Therefore, by Proposition~\ref{PDIpDCH}, $\mathrm{ip}^{(2,X)@}_{s}(\mathrm{CH}^{(2)}_{s}(\mathfrak{Q}^{(2)}))$ is a second-order path in $[\mathfrak{Q}^{(2)}]^{}_{s}$, then, by Lemma~\ref{LDCHNEch}, we have that 
\begin{itemize}
\item[(iii)] $\mathrm{ip}^{(2,X)@}_{s}(\mathrm{CH}^{(2)}_{s}(\mathfrak{Q}^{(2)}))$ is an echelonless second-order path.
\end{itemize}

So, considering the foregoing, by items (i), (ii) and (iii), we have that the second-order path $\mathfrak{R}^{(2)}$, given by
\begin{multline*}
\mathfrak{R}^{(2)}=
\mathrm{ip}^{(2,X)@}_{s}
\left(
\mathrm{CH}^{(2)}_{s}\left(
\mathfrak{Q}^{(2)}
\right)
\right)
\circ^{1\mathbf{Pth}_{\boldsymbol{\mathcal{A}}^{(2)}}}_{s}
\\
\mathrm{ip}^{(2,X)@}_{s}
\left(
\mathrm{CH}^{(2)}_{s}\left(
\mathfrak{P}^{(2),i+1,\bb{\mathfrak{P}^{(2)}}-1}
\right)
\right)
\circ^{1\mathbf{Pth}_{\boldsymbol{\mathcal{A}}^{(2)}}}_{s}
\\
\mathrm{ip}^{(2,X)@}_{s}
\left(
\mathrm{CH}^{(2)}_{s}\left(
\mathfrak{P}^{(2),0,i}
\right)
\right)
\end{multline*}
is an echelonless second-order path that is not head-constant satisfying that $i\in\bb{\mathfrak{P}^{(2)}}-1$ is the greatest index for which $\mathfrak{R}^{(2),0,i}$, which is precisely
$$
\mathfrak{R}^{(2),0,i}=
\mathrm{ip}^{(2,X)@}_{s}\left(
\mathrm{CH}^{(2)}_{s}\left(
\mathfrak{P}^{(2),0,i}
\right)\right),
$$
is a head-constant echelonless second-order path. 

Therefore, we are in position to conclude the initial discussion at the beginning of this subcase. The following chain of equalities holds

\begin{flushleft}
$\mathrm{CH}^{(2)}_{s}
\left(
\mathrm{ip}^{(2,X)@}_{s}
\left(
\mathrm{CH}^{(2)}_{s}\left(
\mathfrak{Q}^{(2)}
\right)
\circ^{1\mathbf{T}_{\Sigma^{\boldsymbol{\mathcal{A}}^{(2)}}}(X)}_{s}
\mathrm{CH}^{(2)}_{s}\left(
\mathfrak{P}^{(2)}
\right)
\right)
\right)
$
\allowdisplaybreaks
\begin{align*}
\quad&=
\mathrm{CH}^{(2)}_{s}
\Bigg(
\mathrm{ip}^{(2,X)@}_{s}
\left(
\mathrm{CH}^{(2)}_{s}\left(
\mathfrak{Q}^{(2)}
\right)
\right)
\circ^{1\mathbf{Pth}_{\boldsymbol{\mathcal{A}}^{(2)}}}_{s}
\\&\qquad\qquad\qquad\qquad
\mathrm{ip}^{(2,X)@}_{s}
\left(
\mathrm{CH}^{(2)}_{s}\left(
\mathfrak{P}^{(2),i+1,\bb{\mathfrak{P}^{(2)}}-1}
\right)
\right)
\circ^{1\mathbf{Pth}_{\boldsymbol{\mathcal{A}}^{(2)}}}_{s}
\\&\qquad\qquad\qquad\qquad\qquad\qquad\qquad\qquad
\mathrm{ip}^{(2,X)@}_{s}
\left(
\mathrm{CH}^{(2)}_{s}\left(
\mathfrak{P}^{(2),0,i}
\right)
\right)
\Bigg)
\tag{1}
\\&=
\mathrm{CH}^{(2)}
\left(
\mathrm{ip}^{(2,X)@}_{s}
\left(
\mathrm{CH}^{(2)}_{s}\left(
\mathfrak{Q}^{(2)}
\right)
\right)
\circ^{1\mathbf{Pth}_{\boldsymbol{\mathcal{A}}^{(2)}}}_{s}
\mathrm{ip}^{(2,X)@}_{s}
\left(
\mathrm{CH}^{(2)}_{s}\left(
\mathfrak{P}^{(2),i,\bb{\mathfrak{P}^{(2)}}-1}
\right)
\right)
\right)
\\&\qquad\qquad\qquad\qquad\qquad\quad
\circ^{1\mathbf{T}_{\Sigma^{\boldsymbol{\mathcal{A}}^{(2)}}}(X)}_{s}
\mathrm{CH}^{(2)}_{s}
\left(
\mathrm{ip}^{(2,X)@}_{s}
\left(
\mathrm{CH}^{(2)}_{s}\left(
\mathfrak{P}^{(2),0,i}
\right)
\right)
\right)
\tag{2}
\\&=
\mathrm{CH}^{(2)}_{s}
\left(
\mathrm{ip}^{(2,X)@}_{s}
\left(
\mathrm{CH}^{(2)}_{s}\left(
\mathfrak{Q}^{(2)}
\right)
\circ^{1\mathbf{T}_{\Sigma^{\boldsymbol{\mathcal{A}}^{(2)}}}(X)}_{s}
\mathrm{CH}^{(2)}_{s}\left(
\mathfrak{P}^{(2),i+1,\bb{\mathfrak{P}^{(2)}}-1}
\right)
\right)
\right)
\\&\qquad\qquad\qquad\qquad\qquad\qquad\qquad\qquad\qquad\quad
\circ^{1\mathbf{T}_{\Sigma^{\boldsymbol{\mathcal{A}}^{(2)}}}(X)}_{s}
\mathrm{CH}^{(2)}_{s}\left(
\mathfrak{P}^{(2),0,i}
\right)
\tag{3}
\\&=
\mathrm{CH}^{(2)}_{s}\left(
\mathfrak{Q}^{(2)}
\circ^{1\mathbf{Pth}_{\boldsymbol{\mathcal{A}}^{(2)}}}_{s}
\mathfrak{P}^{(2),i+1,\bb{\mathfrak{P}^{(2)}}-1}
\right)
\circ^{1\mathbf{T}_{\Sigma^{\boldsymbol{\mathcal{A}}^{(2)}}}(X)}_{s}
\mathrm{CH}^{(2)}_{s}\left(
\mathfrak{P}^{(2),0,i}
\right)
\tag{4}
\\&=
\mathrm{CH}^{(2)}_{s}\left(
\mathfrak{Q}^{(2)}
\circ^{1\mathbf{Pth}_{\boldsymbol{\mathcal{A}}^{(2)}}}_{s}
\mathfrak{P}^{(2)}
\right).
\tag{5}
\end{align*}
\end{flushleft}

The first equality was already proven at the beginning of this subcase;
the second equality simply unravels the second-order Curry-Howard mapping at $\mathfrak{R}^{(2)}$; the third equality follows, on one side from the fact that, by Definition~\ref{DDIp}, $\mathrm{ip}^{(2,X)@}$ is a $\Sigma^{\boldsymbol{\mathcal{A}}^{(2)}}$-homomorphism and, on the other side, from the fact that, by Proposition~\ref{PDIpDCH}, $\mathrm{ip}^{(2,X)@}_{s}(\mathrm{CH}^{(2)}_{s}(\mathfrak{P}^{(2),0,i}))$ is a second-order path in $[\mathfrak{P}^{(2),0,i}]^{}_{s}$; the fourth equality follows from induction. Let us note that the pair $(\mathfrak{Q}^{(2)}\circ^{1\mathbf{Pth}_{\boldsymbol{\mathcal{A}}^{(2)}}}_{s}\mathfrak{P}^{(2),i+1,\bb{\mathfrak{P}^{(2)}}-1},s)$ $\prec_{\mathbf{Pth}_{\boldsymbol{\mathcal{A}}^{(2)}}}$-precedes $(\mathfrak{Q}^{(2)}\circ^{1\mathbf{Pth}_{\boldsymbol{\mathcal{A}}^{(2)}}}_{s}\mathfrak{P}^{(2)},s)$, thus we have that
\begin{multline*}
\mathrm{CH}^{(2)}_{s}
\left(
\mathrm{ip}^{(2,X)@}_{s}
\left(
\mathrm{CH}^{(2)}_{s}\left(
\mathfrak{Q}^{(2)}
\right)
\circ^{1\mathbf{T}_{\Sigma^{\boldsymbol{\mathcal{A}}^{(2)}}}(X)}_{s}
\mathrm{CH}^{(2)}_{s}\left(
\mathfrak{P}^{(2),i+1,\bb{\mathfrak{P}^{(2)}}-1}
\right)
\right)
\right)\\ =
\mathrm{CH}^{(2)}_{s}\left(
\mathfrak{Q}^{(2)}
\circ^{1\mathbf{Pth}_{\boldsymbol{\mathcal{A}}^{(2)}}}_{s}
\mathfrak{P}^{(2),i+1,\bb{\mathfrak{P}^{(2)}}-1}
\right);
\end{multline*}
finally, the last equality recovers the value of the second-order Curry-Howard mapping at $\mathfrak{Q}^{(2)}\circ^{1\mathbf{Pth}_{\boldsymbol{\mathcal{A}}^{(2)}}}_{s}\mathfrak{P}^{(2)}$.

The case $i\in\bb{\mathfrak{P}^{(2)}}-1$ follows.

If~(2.1.2), i.e.,  if $\mathfrak{Q}^{(2)}\circ^{1\mathbf{Pth}_{\boldsymbol{\mathcal{A}}^{(2)}}}_{s}
\mathfrak{P}^{(2)}$ is an echelonless second-order path that is not head-constant and  $i=\bb{\mathfrak{P}^{(2)}}-1$ is the greatest index for which $(\mathfrak{Q}^{(2)}\circ^{1\mathbf{Pth}_{\boldsymbol{\mathcal{A}}^{(2)}}}_{s}\mathfrak{P}^{(2)})^{0,i}$ is head-constant, then, regarding the second-order paths $\mathfrak{P}^{(2)}$ and $\mathfrak{Q}^{(2)}$ we have that
\begin{itemize}
\item[(i)] $\mathfrak{P}^{(2)}$ is a head-constant echelonless second-order path;
\item[(ii)] $\mathfrak{Q}^{(2)}$ is an echelonless second-order path;
\end{itemize}

In this case, according to Definition~\ref{DDCH}, the value of the second-order Curry-Howard mapping at $\mathfrak{Q}^{(2)}\circ^{1\mathbf{Pth}_{\boldsymbol{\mathcal{A}}^{(2)}}}_{s}\mathfrak{P}^{(2)}$ is given by
$$
\mathrm{CH}^{(2)}_{s}\left(
\mathfrak{Q}^{(2)}\circ^{1\mathbf{Pth}_{\boldsymbol{\mathcal{A}}^{(2)}}}_{s}\mathfrak{P}^{(2)}
\right)
=
\mathrm{CH}^{(2)}_{s}\left(
\mathfrak{Q}^{(2)}
\right)
\circ^{1\mathbf{T}_{\Sigma^{\boldsymbol{\mathcal{A}}^{(2)}}}(X)}_{s}
\mathrm{CH}^{(2)}_{s}\left(
\mathfrak{P}^{(2)}
\right).
$$

The following equality is a direct consequence of the equality proved in the first item of this lemma
\begin{multline*}
\mathrm{CH}^{(2)}_{s}
\left(
\mathrm{ip}^{(2,X)@}_{s}
\left(
\mathrm{CH}^{(2)}_{s}\left(
\mathfrak{Q}^{(2)}
\right)
\circ^{1\mathbf{T}_{\Sigma^{\boldsymbol{\mathcal{A}}^{(2)}}}(X)}_{s}
\mathrm{CH}^{(2)}_{s}\left(
\mathfrak{P}^{(2)}
\right)
\right)
\right)
\\=
\mathrm{CH}^{(2)}_{s}
\left(
\mathrm{ip}^{(2,X)@}_{s}
\left(
\mathrm{CH}^{(2)}_{s}\left(
\mathfrak{Q}^{(2)}
\right)
\right)
\circ^{1\mathbf{Pth}_{\boldsymbol{\mathcal{A}}^{(2)}}}_{s}
\mathrm{ip}^{(2,X)@}_{s}
\left(
\mathrm{CH}^{(2)}_{s}\left(
\mathfrak{P}^{(2)}
\right)
\right)
\right)
.
\end{multline*}

We will center our efforts in studying the nature of the second-order path $\mathfrak{R}^{(2)}$ given by
\allowdisplaybreaks
\begin{align*}
\mathfrak{R}^{(2)}&=
\mathrm{ip}^{(2,X)@}_{s}
\left(
\mathrm{CH}^{(2)}_{s}\left(
\mathfrak{Q}^{(2)}
\right)
\right)
\circ^{1\mathbf{Pth}_{\boldsymbol{\mathcal{A}}^{(2)}}}_{s}
\mathrm{ip}^{(2,X)@}_{s}
\left(
\mathrm{CH}^{(2)}_{s}\left(
\mathfrak{P}^{(2)}
\right)
\right).
\end{align*}

Note that
\begin{itemize}
\item[(i)]  $\mathfrak{P}^{(2)}$ is a head-constant echelonless second-order path. 
\end{itemize}

Therefore, for a unique word $\mathbf{s}\in S^{\star}-\{\lambda\}$ and a unique operation symbol $\tau\in\Sigma^{\boldsymbol{\mathcal{A}}}_{\mathbf{s},s}$, the family of first-order translations occurring in $\mathfrak{P}^{(2)}$ is a family of first-order translations of type $\tau$.

By Lemma~\ref{LDCHNEchHd}, we have that $\mathrm{CH}^{(2)}_{s}(\mathfrak{P}^{(2)})\in\mathcal{T}(\tau,\mathrm{T}_{\Sigma^{\boldsymbol{\mathcal{A}}^{(2)}}}(X))^{\star}$, which is a subset of $[\mathrm{T}_{\Sigma^{\boldsymbol{\mathcal{A}}^{(2)}}}(X)]^{\mathsf{HdC}}_{s}$. Therefore, by Proposition~\ref{PDIpDCH}, $\mathrm{ip}^{(2,X)@}_{s}(\mathrm{CH}^{(2)}_{s}(\mathfrak{P}^{(2)}))$ is a second-order path in $[\mathfrak{P}^{(2)}]^{}_{s}$, then, by Lemma~\ref{LDCHNEchHd}, we have that
\begin{itemize}
\item[(i)] $\mathrm{ip}^{(2,X)@}_{s}(\mathrm{CH}^{(2)}_{s}(\mathfrak{P}^{(2)}))$ is a head-constant echelonless second-order path associated to the operation symbol $\tau$.
\end{itemize}

On the other hand 
\begin{itemize}
\item[(ii)] $\mathfrak{Q}^{(2)}$ is an echelonless second-order path.
\end{itemize}

Therefore, for a unique word $\mathbf{s}'\in S^{\star}-\{\lambda\}$ and a unique operation symbol $\tau'\in\Sigma^{\boldsymbol{\mathcal{A}}}_{\mathbf{s}',s}$, the first first-order translation occurring in $\mathfrak{Q}^{(2)}$ is a first-order translation of type $\tau'$.

Let us recall that $\tau\neq\tau'$ since $\mathfrak{Q}^{(2)}\circ^{1\mathbf{Pth}_{\boldsymbol{\mathcal{A}}^{(2)}}}\mathfrak{P}^{(2)}$ is a head-constant echelonless second-order path and $i=\bb{\mathfrak{P}^{(2)}}-1$ is the greatest index for which $(\mathfrak{Q}^{(2)}\circ^{1\mathbf{Pth}_{\boldsymbol{\mathcal{A}}^{(2)}}}\mathfrak{P}^{(2)})^{0,i}$ is a head-constant echelonless second-order path. 

By Lemma~\ref{LDCHNEch}, we have that 
$$\mathrm{CH}^{(2)}_{s}(\mathfrak{Q}^{(2)})\in\mathrm{T}_{\Sigma^{\boldsymbol{\mathcal{A}}^{(2)}}}(X)_{s}
-\left(
\eta^{(2,1)\sharp}\left[
\mathrm{PT}_{\boldsymbol{\mathcal{A}}}
\right]_{s}
\cup
\eta^{(2,\mathcal{A}^{(2)})}\left[\mathcal{A}^{(2)}\right]^{\mathrm{pct}}_{s}
\right).
$$
Therefore, by Proposition~\ref{PDIpDCH}, $\mathrm{ip}^{(2,X)@}_{s}(\mathrm{CH}^{(2)}_{s}(\mathfrak{Q}^{(2)}))$ is a second-order path in $[\mathfrak{Q}^{(2)}]^{}_{s}$, then, by Lemma~\ref{LDCHNEch}, we have that 
\begin{itemize}
\item[(ii)] $\mathrm{ip}^{(2,X)@}_{s}(\mathrm{CH}^{(2)}_{s}(\mathfrak{Q}^{(2)}))$ is an echelonless second-order path whose first first-order translation is associated to the operation symbol $\tau'$.
\end{itemize}

So, considering the foregoing, we have that, by items (i), (ii), the second-order path $\mathfrak{R}^{(2)}$, given by
$$\mathfrak{R}^{(2)}=\mathrm{ip}^{(2,X)@}_{s}\left(
\mathrm{CH}^{(2)}_{s}\left(
\mathfrak{Q}^{(2)}
\right)\right)
\circ^{1\mathbf{Pth}_{\boldsymbol{\mathcal{A}}^{(2)}}}_{s}
\mathrm{ip}^{(2,X)@}_{s}\left(
\mathrm{CH}^{(2)}_{s}\left(
\mathfrak{P}^{(2)}
\right)\right)$$
is an echelonless second-order path that is not head-constant and for which $i=\bb{\mathfrak{P}^{(2)}}-1$ is the greatest index for which $\mathfrak{R}^{(2),0,i}$, which is, precisely,
$$
\mathfrak{R}^{(2),0,i}=
\mathrm{ip}^{(2,X)@}_{s}\left(
\mathrm{CH}^{(2)}_{s}\left(
\mathfrak{P}^{(2)}
\right)\right)
$$
is a head-constant echelonless second-order path. 

Therefore, we are in position to conclude the initial discussion at the beginning of this subcase. The following chain of equalities holds
\begin{flushleft}
$\mathrm{CH}^{(2)}_{s}
\left(
\mathrm{ip}^{(2,X)@}_{s}
\left(
\mathrm{CH}^{(2)}_{s}\left(
\mathfrak{Q}^{(2)}
\right)
\circ^{1\mathbf{T}_{\Sigma^{\boldsymbol{\mathcal{A}}^{(2)}}}(X)}_{s}
\mathrm{CH}^{(2)}_{s}\left(
\mathfrak{P}^{(2)}
\right)
\right)
\right)
$
\allowdisplaybreaks
\begin{align*}
\quad&=
\mathrm{CH}^{(2)}_{s}
\left(
\mathrm{ip}^{(2,X)@}_{s}
\left(
\mathrm{CH}^{(2)}_{s}\left(
\mathfrak{Q}^{(2)}
\right)
\right)
\circ^{1\mathbf{Pth}_{\boldsymbol{\mathcal{A}}^{(2)}}}_{s}
\mathrm{ip}^{(2,X)@}_{s}
\left(
\mathrm{CH}^{(2)}_{s}\left(
\mathfrak{P}^{(2)}
\right)
\right)
\right)
\tag{1}
\\&=
\mathrm{CH}^{(2)}_{s}
\left(
\mathrm{ip}^{(2,X)@}_{s}
\left(
\mathrm{CH}^{(2)}_{s}\left(
\mathfrak{Q}^{(2)}
\right)
\right)
\right)
\circ^{1\mathbf{T}_{\Sigma^{\boldsymbol{\mathcal{A}}^{(2)}}}(X)}_{s}
\\&\qquad\qquad\qquad\qquad\qquad\qquad\qquad\qquad
\mathrm{CH}^{(2)}_{s}
\left(
\mathrm{ip}^{(2,X)@}_{s}
\left(
\mathrm{CH}^{(2)}_{s}\left(
\mathfrak{P}^{(2)}
\right)
\right)
\right)
\tag{2}
\\&=
\mathrm{CH}^{(2)}_{s}\left(
\mathfrak{Q}^{(2)}
\right)
\circ^{1\mathbf{T}_{\Sigma^{\boldsymbol{\mathcal{A}}^{(2)}}}(X)}_{s}
\mathrm{CH}^{(2)}_{s}\left(
\mathfrak{P}^{(2)}
\right)
\tag{3}
\\&=
\mathrm{CH}^{(2)}_{s}\left(
\mathfrak{Q}^{(2)}\circ^{1\mathbf{Pth}_{\boldsymbol{\mathcal{A}}^{(2)}}}_{s}\mathfrak{P}^{(2)}
\right).
\tag{4}
\end{align*}
\end{flushleft}

The first equality was already proven at the beginning of this subcase; the second equality simply unravels the second-order Curry-Howard mapping at the second-order path $\mathfrak{R}^{(2)}$; the third equality follows from the fact that, by Proposition~\ref{PDIpDCH}, the elements $\mathrm{ip}^{(2,X)@}_{s}(
\mathrm{CH}^{(2)}_{s}(
\mathfrak{Q}^{(2)}
))$ and $\mathrm{ip}^{(2,X)@}_{s}(
\mathrm{CH}^{(2)}_{s}(
\mathfrak{P}^{(2)}
))$ are second-order paths in, respectively, $[\mathfrak{Q}^{(2)}]^{}_{s}$ and $[\mathfrak{P}^{(2)}]^{}_{s}$; finally, the last equality recovers the value of the second-order Curry-Howard mapping at $\mathfrak{Q}^{(2)}\circ^{1\mathbf{Pth}_{\boldsymbol{\mathcal{A}}^{(2)}}}_{s}\mathfrak{P}^{(2)}$.

The case $i=\bb{\mathfrak{P}^{(2)}}-1$ follows.

If~(2.1.3), i.e.,  if $\mathfrak{Q}^{(2)}\circ^{1\mathbf{Pth}_{\boldsymbol{\mathcal{A}}^{(2)}}}_{s}
\mathfrak{P}^{(2)}$ is an echelonless second-order path that is not head-constant and $i\in[ \bb{\mathfrak{P}^{(2)}}, 
\bb{\mathfrak{Q}^{(2)}\circ^{1\mathbf{Pth}_{\boldsymbol{\mathcal{A}}^{(2)}}}_{s}\mathfrak{P}^{(2)}}-1
]$ is the greatest index for which $(\mathfrak{Q}^{(2)}\circ^{1\mathbf{Pth}_{\boldsymbol{\mathcal{A}}^{(2)}}}_{s}\mathfrak{P}^{(2)})^{0,i}$ is head-constant, then, regarding the second-order paths $\mathfrak{P}^{(2)}$ and $\mathfrak{Q}^{(2)}$ we have that
\begin{itemize}
\item[(i)] $\mathfrak{P}^{(2)}$ is a head-constant echelonless second-order path;
\item[(ii)] $\mathfrak{Q}^{(2),0,i-\bb{\mathfrak{P}^{(2)}}}$ is a head-constant echelonless second-order path. Moreover, $i-\bb{\mathfrak{P}^{(2)}}\in\bb{\mathfrak{Q}^{(2)}}-1$ is the greatest index for which $\mathfrak{Q}^{(2),0,i-\bb{\mathfrak{P}^{(2)}}}$ is head-constant;
\item[(iii)] $\mathfrak{Q}^{(2),i-\bb{\mathfrak{P}^{(2)}}+1,\bb{\mathfrak{Q}^{(2)}}-1}$ is an echelonless second-order path.
\end{itemize}

This case follows by a similar argument to that presented in Case (2.1.1). We leave the details for the interested reader.

This completes Case~(2.1).

If~(2.2), i.e., if $\mathfrak{Q}^{(2)}\circ^{1\mathbf{Pth}_{\boldsymbol{\mathcal{A}}^{(2)}}}_{s}\mathfrak{P}^{(2)}$ is a head-constant echelonless second-order path that is not coherent, then we let $i\in\bb{\mathfrak{Q}^{(2)}\circ^{1\mathbf{Pth}_{\boldsymbol{\mathcal{A}}^{(2)}}}_{s}\mathfrak{P}^{(2)}}-1$ be the greatest index for which $(\mathfrak{Q}^{(2)}\circ^{1\mathbf{Pth}_{\boldsymbol{\mathcal{A}}^{(2)}}}_{s}\mathfrak{P}^{(2)})^{0,i}$ is a  coherent  head-constant echelonless second-order path. Since $\bb{\mathfrak{Q}^{(2)}\circ^{1\mathbf{Pth}_{\boldsymbol{\mathcal{A}}^{(2)}}}_{s}\mathfrak{P}^{(2)}}=\bb{\mathfrak{Q}^{(2)}}+\bb{\mathfrak{P}^{(2)}}$ then either (2.2.1) $i\in\bb{\mathfrak{P}^{(2)}}-1$, (2.2.2) $i=\bb{\mathfrak{P}^{(2)}}-1$, or (2.2.3) $i\in [
\bb{\mathfrak{P}^{(2)}}, \bb{\mathfrak{Q}^{(2)}
\circ^{1\mathbf{Pth}_{\boldsymbol{\mathcal{A}}^{(2)}}}_{s}
\mathfrak{P}^{(2)}}-1
]$.

If~(2.2.1), i.e., if $\mathfrak{Q}^{(2)}\circ^{1\mathbf{Pth}_{\boldsymbol{\mathcal{A}}^{(2)}}}_{s}
\mathfrak{P}^{(2)}$ is a head-constant echelonless second-order path that is not coherent and $i\in\bb{\mathfrak{P}^{(2)}}$ is the greatest index for which $(\mathfrak{Q}^{(2)}\circ^{1\mathbf{Pth}_{\boldsymbol{\mathcal{A}}^{(2)}}}_{s}
\mathfrak{P}^{(2)})^{0,i}$ is a coherent head-constant echelonless second-order path then, regarding the second-order paths $\mathfrak{Q}^{(2)}$ and $\mathfrak{P}^{(2)}$, we have that
\begin{itemize}
\item[(i)] $\mathfrak{P}^{(2),0,i}$ is a  coherent head-constant echelonless second-order path. Moreover, $i\in\bb{\mathfrak{P}^{(2)}}-1$ is the greatest index for which $\mathfrak{P}^{(2),0,i}$ is coherent;
\item[(ii)] $\mathfrak{P}^{(2),i+1,\bb{\mathfrak{P}^{(2)}}-1}$ is a head-constant echelonless second-order path;
\item[(iii)] $\mathfrak{Q}^{(2)}$ is a head-constant echelonless second-order path.
\end{itemize}

Following (i) and (ii) and taking into account Definition~\ref{DDCH}, we have that the value of the second-order Curry-Howard mapping at $\mathfrak{P}^{(2)}$ is given by
$$
\mathrm{CH}^{(2)}_{s}\left(
\mathfrak{P}^{(2)}
\right)
=
\mathrm{CH}^{(2)}_{s}\left(
\mathfrak{P}^{(2),i+1,\bb{\mathfrak{P}^{(2)}}-1}
\right)
\circ^{1\mathbf{T}_{\Sigma^{\boldsymbol{\mathcal{A}}^{(2)}}}(X)}_{s}
\mathrm{CH}^{(2)}_{s}\left(
\mathfrak{P}^{(2),0,i}
\right).
$$

Moreover, we have that 
$$
\left(\mathfrak{Q}^{(2)}
\circ^{1\mathbf{Pth}_{\boldsymbol{\mathcal{A}}^{(2)}}}_{s}
\mathfrak{P}^{(2)}
\right)^{i+1,\bb{\mathfrak{Q}^{(2)}\circ^{1\mathbf{Pth}_{\boldsymbol{\mathcal{A}}^{(2)}}}_{s}\mathfrak{P}^{(2)}}-1}
=
\mathfrak{Q}^{(2)}
\circ^{1\mathbf{Pth}_{\boldsymbol{\mathcal{A}}^{(2)}}}_{s}
\mathfrak{P}^{(2),i+1,\bb{\mathfrak{P}^{(2)}}-1}.
$$

Thus, the value of the second-order Curry-Howard mapping at $\mathfrak{Q}^{(2)}
\circ^{1\mathbf{Pth}_{\boldsymbol{\mathcal{A}}^{(2)}}}_{s}
\mathfrak{P}^{(2)}$ is given by
\begin{multline*}
\mathrm{CH}^{(2)}_{s}\left(
\mathfrak{Q}^{(2)}
\circ^{1\mathbf{Pth}_{\boldsymbol{\mathcal{A}}^{(2)}}}_{s}
\mathfrak{P}^{(2)}
\right)
\\=
\mathrm{CH}^{(2)}_{s}\left(
\mathfrak{Q}^{(2)}
\circ^{1\mathbf{Pth}_{\boldsymbol{\mathcal{A}}^{(2)}}}_{s}
\mathfrak{P}^{(2),i+1,\bb{\mathfrak{P}^{(2)}}-1}
\right)
\circ^{1\mathbf{T}_{\Sigma^{\boldsymbol{\mathcal{A}}^{(2)}}}(X)}_{s}
\mathrm{CH}^{(2)}\left(
\mathfrak{P}^{(2),0,i}
\right).
\end{multline*}

On the other hand, the following chain of equalities holds
\begin{flushleft}
$\mathrm{CH}^{(2)}_{s}
\left(
\mathrm{ip}^{(2,X)@}_{s}
\left(
\mathrm{CH}^{(2)}_{s}\left(
\mathfrak{Q}^{(2)}
\right)
\circ^{1\mathbf{T}_{\Sigma^{\boldsymbol{\mathcal{A}}^{(2)}}}(X)}_{s}
\mathrm{CH}^{(2)}_{s}\left(
\mathfrak{P}^{(2)}
\right)
\right)
\right)
$
\allowdisplaybreaks
\begin{align*}
\quad&=
\mathrm{CH}^{(2)}_{s}
\left(
\mathrm{ip}^{(2,X)@}_{s}
\left(
\mathrm{CH}^{(2)}_{s}\left(
\mathfrak{Q}^{(2)}
\right)
\right)
\circ^{1\mathbf{Pth}_{\boldsymbol{\mathcal{A}}^{(2)}}}_{s}
\mathrm{ip}^{(2,X)@}_{s}
\left(
\mathrm{CH}^{(2)}_{s}\left(
\mathfrak{P}^{(2)}
\right)
\right)
\right)
\tag{1}
\\&=
\mathrm{CH}^{(2)}_{s}
\Bigg(
\mathrm{ip}^{(2,X)@}_{s}
\left(
\mathrm{CH}^{(2)}_{s}\left(
\mathfrak{Q}^{(2)}
\right)
\right)
\circ^{1\mathbf{Pth}_{\boldsymbol{\mathcal{A}}^{(2)}}}_{s}
\\&\qquad\quad
\mathrm{ip}^{(2,X)@}_{s}
\left(
\mathrm{CH}^{(2)}_{s}\left(
\mathfrak{P}^{(2),i+1,\bb{\mathfrak{P}^{(2)}}-1}
\right)
\circ^{1\mathbf{T}_{\Sigma^{\boldsymbol{\mathcal{A}}^{(2)}}}(X)}_{s}
\mathrm{CH}^{(2)}_{s}\left(
\mathfrak{P}^{(2),0,i}
\right)
\right)
\Bigg)
\tag{2}
\\&=
\mathrm{CH}^{(2)}_{s}
\Bigg(
\mathrm{ip}^{(2,X)@}_{s}
\left(
\mathrm{CH}^{(2)}_{s}\left(
\mathfrak{Q}^{(2)}
\right)
\right)
\circ^{1\mathbf{Pth}_{\boldsymbol{\mathcal{A}}^{(2)}}}_{s}
\\&\qquad\qquad\qquad\qquad
\mathrm{ip}^{(2,X)@}_{s}
\left(
\mathrm{CH}^{(2)}_{s}\left(
\mathfrak{P}^{(2),i+1,\bb{\mathfrak{P}^{(2)}}-1}
\right)
\right)
\circ^{1\mathbf{Pth}_{\boldsymbol{\mathcal{A}}^{(2)}}}_{s}
\\&\qquad\qquad\qquad\qquad\qquad\qquad\qquad\qquad
\mathrm{ip}^{(2,X)@}_{s}
\left(
\mathrm{CH}^{(2)}_{s}\left(
\mathfrak{P}^{(2),0,i}
\right)
\right)
\Bigg).
\tag{3}
\end{align*}
\end{flushleft}

The first equality is a direct consequence of the equation proven in the first item of this lemma; the second equality unravels the value of the second-order Curry-Howard mapping at $\mathfrak{P}^{(2)}$; the third equality follows from the fact that, by Definition~\ref{DDIp}, $\mathrm{ip}^{(2,X)@}$ is a $\Sigma^{\boldsymbol{\mathcal{A}}^{(2)}}$-homomorphism and from the fact that, by Proposition~\ref{PDIpDCH}, the elements $\mathrm{ip}^{(2,X)@}_{s}(\mathrm{CH}^{(2)}_{s}(\mathfrak{P}^{(2),i,\bb{\mathfrak{P}^{(2)}}-1}))$ and $\mathrm{ip}^{(2,X)@}_{s}(\mathrm{CH}^{(2)}_{s}(\mathfrak{P}^{(2),0,i}))$ are second-order paths in $\mathrm{Pth}_{\boldsymbol{\mathcal{A}}^{(2)},s}$. We have omitted parenthesis because, by Proposition~\ref{PDPthComp}, the $1$-composition of second-order paths is associative.

We will center our efforts in studying the nature of the second-order path $\mathfrak{R}^{(2)}$, given by
\begin{multline*}
\mathfrak{R}^{(2)}=
\mathrm{ip}^{(2,X)@}_{s}\left(
\mathrm{CH}^{(2)}_{s}\left(
\mathfrak{Q}^{(2)}
\right)\right)
\circ^{1\mathbf{Pth}_{\boldsymbol{\mathcal{A}}^{(2)}}}_{s}
\\
\mathrm{ip}^{(2,X)@}_{s}\left(
\mathrm{CH}^{(2)}_{s}\left(
\mathfrak{P}^{(2),i+1,\bb{\mathfrak{P}^{(2)}}-1}
\right)\right)
\\
\circ^{1\mathbf{Pth}_{\boldsymbol{\mathcal{A}}^{(2)}}}_{s}
\mathrm{ip}^{(2,X)@}_{s}\left(
\mathrm{CH}^{(2)}_{s}\left(
\mathfrak{P}^{(2),0,i}
\right)\right).
\end{multline*}

Note that
\begin{itemize}
\item[(i)]  $\mathfrak{P}^{(2)}$ is a head-constant echelonless second-order path and $i\in\bb{\mathfrak{P}^{(2)}}-1$ is the greatest index for which $\mathfrak{P}^{(2),0,i}$ is a coherent echelonless second-order path.
\end{itemize}

Moreover, in virtue of Proposition~\ref{PDIpDCH}, $\mathrm{ip}^{(2,X)@}_{s}(\mathrm{CH}^{(2)}_{s}(\mathfrak{P}^{(2),0,i}))$ and  $\mathrm{ip}^{(2,X)@}_{s}(\mathrm{CH}^{(2)}_{s}(\mathfrak{P}^{(2),i+1,\bb{\mathfrak{P}^{(2)}}-1}))$  are second-order paths in $[\mathfrak{P}^{(2),0,i}]^{}_{s}$ and $[\mathfrak{P}^{(2),i+1,\bb{\mathfrak{P}^{(2)}}-1}]^{}_{s}$, respectively.

In view of the foregoing and taking into account Lemma~\ref{LTech}, we have that 
\begin{itemize}
\item[(i)] $\mathrm{ip}^{(2,X)@}_{s}(\mathrm{CH}^{(2)}_{s}(\mathfrak{P}^{(2),i+1,\bb{\mathfrak{P}^{(2)}}-1}))\circ_{s}^{1\mathbf{Pth}_{\boldsymbol{\mathcal{A}}^{(2)}}}\mathrm{ip}^{(2,X)@}_{s}(\mathrm{CH}^{(2)}_{s}(\mathfrak{P}^{(2),0,i}))$ is a head-constant echelonless second-order path and $i\in\bb{\mathfrak{P}^{(2)}}-1$ is the greatest index for which $(\mathrm{ip}^{(2,X)@}_{s}(\mathrm{CH}^{(2)}_{s}(\mathfrak{P}^{(2),i+1,\bb{\mathfrak{P}^{(2)}}-1}))\circ_{s}^{1\mathbf{Pth}_{\boldsymbol{\mathcal{A}}^{(2)}}}\mathrm{ip}^{(2,X)@}_{s}(\mathrm{CH}^{(2)}_{s}(\mathfrak{P}^{(2),0,i})))^{0,i}$ is a coherent echelonless second-order path.
\end{itemize}

Note also that 
\begin{multline*}
\left(\mathrm{ip}^{(2,X)@}_{s}\left(
\mathrm{CH}^{(2)}_{s}\left(
\mathfrak{P}^{(2),i+1,\bb{\mathfrak{P}^{(2)}}-1}
\right)\right)\circ_{s}^{1\mathbf{Pth}_{\boldsymbol{\mathcal{A}}^{(2)}}}
\mathrm{ip}^{(2,X)@}_{s}\left(
\mathrm{CH}^{(2)}_{s}\left(
\mathfrak{P}^{(2),0,i}
\right)\right)
\right)^{0,i}
\\=
\mathrm{ip}^{(2,X)@}_{s}\left(
\mathrm{CH}^{(2)}_{s}\left(
\mathfrak{P}^{(2),0,i}
\right)\right).
\end{multline*}

Therefore, the following chain of equalities holds
\begin{flushleft}
$\mathrm{CH}^{(2)}_{s}\left(
\mathrm{ip}^{(2,X)@}_{s}\left(
\mathrm{CH}^{(2)}_{s}\left(
\mathfrak{P}^{(2),i+1,\bb{\mathfrak{P}^{(2)}}-1}
\right)\right)
\circ_{s}^{1\mathbf{Pth}_{\boldsymbol{\mathcal{A}}^{(2)}}}
\mathrm{ip}^{(2,X)@}_{s}\left(
\mathrm{CH}^{(2)}_{s}\left(
\mathfrak{P}^{(2),0,i}
\right)\right)
\right)$
\allowdisplaybreaks
\begin{align*}
\quad&=
\mathrm{CH}^{(2)}_{s}\left(
\mathrm{ip}^{(2,X)@}_{s}\left(
\mathrm{CH}^{(2)}_{s}\left(
\mathfrak{P}^{(2),i+1,\bb{\mathfrak{P}^{(2)}}-1}
\right)\right)\right)
\circ_{s}^{1\mathbf{T}_{\Sigma^{\boldsymbol{\mathcal{A}}^{(2)}}}(X)}
\\&\qquad\qquad\qquad\qquad\qquad\qquad\qquad\qquad
\mathrm{CH}^{(2)}_{s}\left(
\mathrm{ip}^{(2,X)@}_{s}\left(
\mathrm{CH}^{(2)}_{s}\left(
\mathfrak{P}^{(2),0,i}
\right)\right)
\right)
\tag{1}
\\&=
\mathrm{CH}^{(2)}_{s}\left(
\mathfrak{P}^{(2),i+1,\bb{\mathfrak{P}^{(2)}}-1}
\right)
\circ_{s}^{1\mathbf{T}_{\Sigma^{\boldsymbol{\mathcal{A}}^{(2)}}}(X)}
\mathrm{CH}^{(2)}_{s}\left(
\mathfrak{P}^{(2),0,i}
\right)
\tag{2}
\\&=
\mathrm{CH}^{(2)}_{s}\left(
\mathfrak{P}^{(2)}
\right).
\tag{3}
\end{align*}
\end{flushleft}

The first equality simply unravels the second-order Curry-Howard mapping at the second-order path under consideration. In this regard we have taken into consideration the previous descriptions;  the second equality follows from Proposition~\ref{PDIpDCH}; finally, the last equality recovers the second-order Curry-Howard mapping at the second-order path $\mathfrak{P}^{(2)}$.

By Lemma~\ref{LDCHNEchHdNC}, we have that $\mathrm{CH}^{(2)}_{s}(\mathfrak{P}^{(2)})$ belongs to $\mathrm{T}_{\Sigma^{\boldsymbol{\mathcal{A}}^{(2)}}}(X)^{\mathsf{HdC}\And\mathsf{C}}_{s}$. Therefore, there exists a unique word $\mathbf{s}\in S^{\star}-\{\lambda\}$ and a unique operation symbol $\tau\in\Sigma^{\boldsymbol{\mathcal{A}}}_{\mathbf{s},s}$ associated to $\mathfrak{P}^{(2)}$. In virtue of the just stated chain of equalities,  the head-constant echelonless second-order path $\mathrm{ip}^{(2,X)@}_{s}(\mathrm{CH}^{(2)}_{s}(\mathfrak{P}^{(2),i+1,\bb{\mathfrak{P}^{(2)}}-1}))\circ_{s}^{1\mathbf{Pth}_{\boldsymbol{\mathcal{A}}^{(2)}}}\mathrm{ip}^{(2,X)@}_{s}(\mathrm{CH}^{(2)}_{s}(\mathfrak{P}^{(2),0,i}))$ is also associated to the same operation symbol $\tau\in\Sigma^{\boldsymbol{\mathcal{A}}}_{\mathbf{s},s}$.

Let us also recall that
\begin{itemize}
\item[(ii)]  $\mathfrak{Q}^{(2)}$ is a head-constant echelonless second-order path. 
\end{itemize}

Therefore, for the unique word $\mathbf{s}\in S^{\star}-\{\lambda\}$ and the unique operation symbol $\tau\in\Sigma^{\boldsymbol{\mathcal{A}}}_{\mathbf{s},s}$ the first first-order translation occurring in $\mathfrak{Q}^{(2)}$ is a first-order translation of type $\tau$.

Note that the operation symbol $\tau$ is the same as in case (i), since $\mathfrak{Q}^{(2)}\circ^{1\mathbf{Pth}_{\boldsymbol{\mathcal{A}}^{(2)}}}_{s}\mathfrak{P}^{(2)}$ is  head-constant by hypothesis.

By Lemma~\ref{LDCHNEchHd}, we have that $\mathrm{CH}^{(2)}_{s}(\mathfrak{Q}^{(2)})\in\mathcal{T}(\tau,\mathrm{T}_{\Sigma^{\boldsymbol{\mathcal{A}}^{(2)}}}(X))^{\star}$, which is a subset of $\mathrm{T}_{\Sigma^{\boldsymbol{\mathcal{A}}^{(2)}}}(X)^{\mathsf{HdC}}_{s}$. By Proposition~\ref{PDIpDCH}, $\mathrm{ip}^{(2,X)@}_{s}(\mathrm{CH}^{(2)}_{s}(\mathfrak{Q}^{(2)}))$ is a second-order path in $[\mathfrak{Q}^{(2)}]^{}_{s}$ then, by Lemma~\ref{LDCHNEchHd}, we have that
\begin{itemize}
\item[(ii)] $\mathrm{ip}^{(2,X)@}_{s}(\mathrm{CH}^{(2)}_{s}(\mathfrak{Q}^{(2)}))$ is a head-constant echelonless second-order path associated to the operation symbol $\tau$.
\end{itemize}

Taking into account items (i) and (ii), we have that the 
the second-order path $\mathfrak{R}^{(2)}$, given by
\begin{multline*}
\mathfrak{R}^{(2)}=
\mathrm{ip}^{(2,X)@}_{s}\left(
\mathrm{CH}^{(2)}_{s}\left(
\mathfrak{Q}^{(2)}
\right)\right)
\circ^{1\mathbf{Pth}_{\boldsymbol{\mathcal{A}}^{(2)}}}_{s}
\\
\mathrm{ip}^{(2,X)@}_{s}\left(
\mathrm{CH}^{(2)}_{s}\left(
\mathfrak{P}^{(2),i+1,\bb{\mathfrak{P}^{(2)}}-1}
\right)\right)
\\
\circ^{1\mathbf{Pth}_{\boldsymbol{\mathcal{A}}^{(2)}}}_{s}
\mathrm{ip}^{(2,X)@}_{s}\left(
\mathrm{CH}^{(2)}_{s}\left(
\mathfrak{P}^{(2),0,i}
\right)\right).
\end{multline*}
is a head-constant echelonless second-order path that is not coherent and for which $i\in\bb{\mathfrak{P}^{(2)}}-1$ is the greatest index for which $\mathfrak{R}^{(2),0,i}$ is coherent, which is precisely
$$
\mathfrak{R}^{(2),0,i}=
\mathrm{ip}^{(2,X)@}_{s}\left(
\mathrm{CH}^{(2)}_{s}\left(
\mathfrak{P}^{(2),0,i}
\right)\right),
$$ is coherent.

Therefore, we are in position to conclude the initial discussion at the beginning of this subcase. The following chain of equalities holds
\begin{flushleft}
$\mathrm{CH}^{(2)}_{s}
\left(
\mathrm{ip}^{(2,X)@}_{s}
\left(
\mathrm{CH}^{(2)}_{s}\left(
\mathfrak{Q}^{(2)}
\right)
\circ^{1\mathbf{T}_{\Sigma^{\boldsymbol{\mathcal{A}}^{(2)}}}(X)}_{s}
\mathrm{CH}^{(2)}_{s}\left(
\mathfrak{P}^{(2)}
\right)
\right)
\right)
$
\allowdisplaybreaks
\begin{align*}
\qquad&=
\mathrm{CH}^{(2)}_{s}
\Bigg(
\mathrm{ip}^{(2,X)@}_{s}
\left(
\mathrm{CH}^{(2)}_{s}\left(
\mathfrak{Q}^{(2)}
\right)
\right)
\circ^{1\mathbf{Pth}_{\boldsymbol{\mathcal{A}}^{(2)}}}_{s}
\\&\qquad\qquad\qquad\qquad
\mathrm{ip}^{(2,X)@}_{s}
\left(
\mathrm{CH}^{(2)}_{s}\left(
\mathfrak{P}^{(2),i+1,\bb{\mathfrak{P}^{(2)}}-1}
\right)
\right)
\circ^{1\mathbf{Pth}_{\boldsymbol{\mathcal{A}}^{(2)}}}_{s}
\\&\qquad\qquad\qquad\qquad\qquad\qquad\qquad\qquad
\mathrm{ip}^{(2,X)@}_{s}
\left(
\mathrm{CH}^{(2)}_{s}\left(
\mathfrak{P}^{(2),0,i}
\right)
\right)
\Bigg)
\tag{1}
\\&=
\mathrm{CH}^{(2)}
\left(
\mathrm{ip}^{(2,X)@}_{s}
\left(
\mathrm{CH}^{(2)}_{s}\left(
\mathfrak{Q}^{(2)}
\right)
\right)
\circ^{1\mathbf{Pth}_{\boldsymbol{\mathcal{A}}^{(2)}}}_{s}
\right.
\\&\qquad\qquad\qquad
\left.
\mathrm{ip}^{(2,X)@}_{s}
\left(
\mathrm{CH}^{(2)}_{s}\left(
\mathfrak{P}^{(2),i+1,\bb{\mathfrak{P}^{(2)}}-1}
\right)
\right)
\right)
\\&\qquad\qquad\qquad\qquad\qquad
\circ^{1\mathbf{T}_{\Sigma^{\boldsymbol{\mathcal{A}}^{(2)}}}(X)}_{s}
\mathrm{CH}^{(2)}_{s}
\left(
\mathrm{ip}^{(2,X)@}_{s}
\left(
\mathrm{CH}^{(2)}_{s}\left(
\mathfrak{P}^{(2),0,i}
\right)
\right)
\right)
\tag{2}
\\&=
\mathrm{CH}^{(2)}_{s}
\left(
\mathrm{ip}^{(2,X)@}_{s}
\left(
\mathrm{CH}^{(2)}_{s}\left(
\mathfrak{Q}^{(2)}
\right)
\circ^{1\mathbf{T}_{\Sigma^{\boldsymbol{\mathcal{A}}^{(2)}}}(X)}_{s}
\mathrm{CH}^{(2)}_{s}\left(
\mathfrak{P}^{(2),i+1,\bb{\mathfrak{P}^{(2)}}-1}
\right)
\right)
\right)
\\&\qquad\qquad\qquad\qquad\qquad\qquad\qquad\qquad\qquad
\circ^{1\mathbf{T}_{\Sigma^{\boldsymbol{\mathcal{A}}^{(2)}}}(X)}_{s}
\mathrm{CH}^{(2)}_{s}\left(
\mathfrak{P}^{(2),0,i}
\right)
\tag{3}
\\&=
\mathrm{CH}^{(2)}_{s}\left(
\mathfrak{Q}^{(2)}
\circ^{1\mathbf{Pth}_{\boldsymbol{\mathcal{A}}^{(2)}}}_{s}
\mathfrak{P}^{(2),i+1,\bb{\mathfrak{P}^{(2)}}-1}
\right)
\circ^{1\mathbf{T}_{\Sigma^{\boldsymbol{\mathcal{A}}^{(2)}}}(X)}_{s}
\mathrm{CH}^{(2)}_{s}\left(
\mathfrak{P}^{(2),0,i}
\right)
\tag{4}
\\&=
\mathrm{CH}^{(2)}_{s}\left(
\mathfrak{Q}^{(2)}
\circ^{1\mathbf{Pth}_{\boldsymbol{\mathcal{A}}^{(2)}}}_{s}
\mathfrak{P}^{(2)}
\right).
\tag{5}
\end{align*}
\end{flushleft}

The first equality was already proven at the beginning of this case; the second equality simply unravels the second-order Curry-Howard mapping at the second-order path $\mathfrak{R}^{(2)}$; the third equality follows, on the one hand, from the fact that, in virtue of Definition~\ref{DDIp}, $\mathrm{ip}^{(2,X)@}$ is a $\Sigma^{\boldsymbol{\mathcal{A}}^{(2)}}$-homomorphism and, on the other hand, from the fact that, by Proposition~\ref{PDIpDCH}, $\mathrm{ip}^{(2,X)@}_{s}(\mathrm{CH}^{(2)}_{s}(\mathfrak{P}^{(2),0,i}))$ is a second-order path in $[\mathfrak{P}^{(2),0,i}]^{}_{s}$; the fourth equality follows from induction. Let us note that the pair $(\mathfrak{Q}^{(2)}\circ^{1\mathbf{Pth}_{\boldsymbol{\mathcal{A}}^{(2)}}}_{s}\mathfrak{P}^{(2),i+1,\bb{\mathfrak{P}^{(2)}}-1},s)$ $\prec_{\mathbf{Pth}_{\boldsymbol{\mathcal{A}}^{(2)}}}$-precedes $(\mathfrak{Q}^{(2)}\circ^{1\mathbf{Pth}_{\boldsymbol{\mathcal{A}}^{(2)}}}_{s}\mathfrak{P}^{(2)},s)$, thus we have that
\begin{multline*}
\mathrm{CH}^{(2)}_{s}
\left(
\mathrm{ip}^{(2,X)@}_{s}
\left(
\mathrm{CH}^{(2)}_{s}\left(
\mathfrak{Q}^{(2)}
\right)
\circ^{1\mathbf{T}_{\Sigma^{\boldsymbol{\mathcal{A}}^{(2)}}}(X)}_{s}
\mathrm{CH}^{(2)}_{s}\left(
\mathfrak{P}^{(2),i+1,\bb{\mathfrak{P}^{(2)}}-1}
\right)
\right)
\right)
\\ =
\mathrm{CH}^{(2)}_{s}\left(
\mathfrak{Q}^{(2)}
\circ^{1\mathbf{Pth}_{\boldsymbol{\mathcal{A}}^{(2)}}}_{s}
\mathfrak{P}^{(2),i+1,\bb{\mathfrak{P}^{(2)}}-1}
\right);
\end{multline*}
finally, the last equality recovers the value of the second-order Curry-Howard mapping at $\mathfrak{Q}^{(2)}\circ^{1\mathbf{Pth}_{\boldsymbol{\mathcal{A}}^{(2)}}}_{s}\mathfrak{P}^{(2)}$.

The case $i\in\bb{\mathfrak{P}^{(2)}}-1$ follows.

If~(2.2.2), i.e., if $\mathfrak{Q}^{(2)}\circ^{1\mathbf{Pth}_{\boldsymbol{\mathcal{A}}^{(2)}}}_{s}
\mathfrak{P}^{(2)}$ is a head-constant echelonless second-order path that is not coherent and $i=\bb{\mathfrak{P}^{(2)}}-1$ is the greatest index for which $(\mathfrak{Q}^{(2)}\circ^{1\mathbf{Pth}_{\boldsymbol{\mathcal{A}}^{(2)}}}_{s}\mathfrak{P}^{(2)})^{0,i}$ is coherent, the value of the second-order Curry-Howard mapping at $\mathfrak{Q}^{(2)}\circ^{1\mathbf{Pth}_{\boldsymbol{\mathcal{A}}^{(2)}}}_{s}\mathfrak{P}^{(2)}$ is given by
$$
\mathrm{CH}^{(2)}_{s}\left(
\mathfrak{Q}^{(2)}\circ^{1\mathbf{Pth}_{\boldsymbol{\mathcal{A}}^{(2)}}}_{s}\mathfrak{P}^{(2)}
\right)
=
\mathrm{CH}^{(2)}_{s}\left(
\mathfrak{Q}^{(2)}
\right)
\circ^{1\mathbf{T}_{\Sigma^{\boldsymbol{\mathcal{A}}^{(2)}}}(X)}_{s}
\mathrm{CH}^{(2)}_{s}\left(
\mathfrak{P}^{(2)}
\right).
$$

The following equality is a direct consequence of the equation proven in the first item of this lemma
\begin{multline*}
\mathrm{CH}^{(2)}_{s}
\left(
\mathrm{ip}^{(2,X)@}_{s}
\left(
\mathrm{CH}^{(2)}_{s}\left(
\mathfrak{Q}^{(2)}
\right)
\circ^{1\mathbf{T}_{\Sigma^{\boldsymbol{\mathcal{A}}^{(2)}}}(X)}_{s}
\mathrm{CH}^{(2)}_{s}\left(
\mathfrak{P}^{(2)}
\right)
\right)
\right)
\\=
\mathrm{CH}^{(2)}_{s}
\left(
\mathrm{ip}^{(2,X)@}_{s}
\left(
\mathrm{CH}^{(2)}_{s}\left(
\mathfrak{Q}^{(2)}
\right)
\right)
\circ^{1\mathbf{Pth}_{\boldsymbol{\mathcal{A}}^{(2)}}}_{s}
\mathrm{ip}^{(2,X)@}_{s}
\left(
\mathrm{CH}^{(2)}_{s}\left(
\mathfrak{P}^{(2)}
\right)
\right)
\right)
.
\end{multline*}

We will study the nature of the second-order path $\mathfrak{R}^{(2)}$ given by
$$
\mathfrak{R}^{(2)}=
\mathrm{ip}^{(2,X)@}_{s}
\left(
\mathrm{CH}^{(2)}_{s}\left(
\mathfrak{Q}^{(2)}
\right)
\right)
\circ^{1\mathbf{Pth}_{\boldsymbol{\mathcal{A}}^{(2)}}}_{s}
\mathrm{ip}^{(2,X)@}_{s}
\left(
\mathrm{CH}^{(2)}_{s}\left(
\mathfrak{P}^{(2)}
\right)
\right).
$$

Note that
\begin{itemize}
\item[(i)]  $\mathfrak{Q}^{(2)}\circ^{1\mathbf{Pth}_{\boldsymbol{\mathcal{A}}^{(2)}}}_{s}\mathfrak{P}^{(2)}$ is a head-constant echelonless second-order path and $i=\bb{\mathfrak{P}^{(2)}}-1$ is the greatest index for which $(\mathfrak{Q}^{(2)}\circ^{1\mathbf{Pth}_{\boldsymbol{\mathcal{A}}^{(2)}}}_{s}\mathfrak{P}^{(2)})^{0,i}$ is a coherent echelonless second-order path.
\end{itemize}

Moreover, in virtue of Proposition~\ref{PDIpDCH}, $\mathrm{ip}^{(2,X)@}_{s}(\mathrm{CH}^{(2)}_{s}(\mathfrak{P}^{(2)}))$ and  $\mathrm{ip}^{(2,X)@}_{s}(\mathrm{CH}^{(2)}_{s}(\mathfrak{Q}^{(2)}))$  are second-order paths in $[\mathfrak{P}^{(2)}]^{}_{s}$ and $[\mathfrak{Q}^{(2)}]^{}_{s}$, respectively.

In view of the foregoing and taking into account Lemma~\ref{LTech}, we have that 
\begin{itemize}
\item[(i)]  $\mathfrak{R}^{(2)}$ is a head-constant echelonless second-order path and $i=\bb{\mathfrak{P}^{(2)}}-1$ is the greatest index for which $\mathfrak{R}^{(2),0,i}$ is a coherent echelonless second-order path. Moreover,
$$
\mathfrak{R}^{(2),0,i}=
\mathrm{ip}^{(2,X)@}_{s}\left(
\mathrm{CH}^{(2)}_{s}\left(
\mathfrak{P}^{(2)}
\right)\right).
$$ 
\end{itemize}

Therefore, we are in position to conclude the initial discussion at the beginning of this subcase. The following chain of equalities holds
\begin{flushleft}
$\mathrm{CH}^{(2)}_{s}
\left(
\mathrm{ip}^{(2,X)@}_{s}
\left(
\mathrm{CH}^{(2)}_{s}\left(
\mathfrak{Q}^{(2)}
\right)
\circ^{1\mathbf{T}_{\Sigma^{\boldsymbol{\mathcal{A}}^{(2)}}}(X)}_{s}
\mathrm{CH}^{(2)}_{s}\left(
\mathfrak{P}^{(2)}
\right)
\right)
\right)
$
\allowdisplaybreaks
\begin{align*}
\quad&=
\mathrm{CH}^{(2)}_{s}
\left(
\mathrm{ip}^{(2,X)@}_{s}
\left(
\mathrm{CH}^{(2)}_{s}\left(
\mathfrak{Q}^{(2)}
\right)
\right)
\circ^{1\mathbf{Pth}_{\boldsymbol{\mathcal{A}}^{(2)}}}_{s}
\mathrm{ip}^{(2,X)@}_{s}
\left(
\mathrm{CH}^{(2)}_{s}\left(
\mathfrak{P}^{(2)}
\right)
\right)
\right)
\tag{1}
\\&=
\mathrm{CH}^{(2)}_{s}
\left(
\mathrm{ip}^{(2,X)@}_{s}
\left(
\mathrm{CH}^{(2)}_{s}\left(
\mathfrak{Q}^{(2)}
\right)
\right)
\right)
\circ^{1\mathbf{T}_{\Sigma^{\boldsymbol{\mathcal{A}}^{(2)}}}(X)}_{s}
\\&\qquad\qquad\qquad\qquad\qquad\qquad\qquad\qquad
\mathrm{CH}^{(2)}_{s}
\left(
\mathrm{ip}^{(2,X)@}_{s}
\left(
\mathrm{CH}^{(2)}_{s}\left(
\mathfrak{P}^{(2)}
\right)
\right)
\right)
\tag{2}
\\&=
\mathrm{CH}^{(2)}_{s}\left(
\mathfrak{Q}^{(2)}
\right)
\circ^{1\mathbf{T}_{\Sigma^{\boldsymbol{\mathcal{A}}^{(2)}}}(X)}_{s}
\mathrm{CH}^{(2)}_{s}\left(
\mathfrak{P}^{(2)}
\right)
\tag{3}
\\&=
\mathrm{CH}^{(2)}_{s}\left(
\mathfrak{Q}^{(2)}\circ^{1\mathbf{Pth}_{\boldsymbol{\mathcal{A}}^{(2)}}}_{s}\mathfrak{P}^{(2)}
\right).
\tag{4}
\end{align*}
\end{flushleft}

The first equality was already proved at the beginning of this subcase; the second equality simply unravels the second-order Curry-Howard mapping at the second-order path $\mathfrak{R}^{(2)}$; the third equality follows from the fact that, by Proposition~\ref{PDIpDCH}, the elements $\mathrm{ip}^{(2,X)@}_{s}(
\mathrm{CH}^{(2)}_{s}(
\mathfrak{Q}^{(2)}
))$ and $\mathrm{ip}^{(2,X)@}_{s}(
\mathrm{CH}^{(2)}_{s}(
\mathfrak{P}^{(2)}
))$ are second-order paths in, respectively, $[\mathfrak{Q}^{(2)}]^{}_{s}$ and $[\mathfrak{P}^{(2)}]^{}_{s}$; finally, the last equality recovers the value of the second-order Curry-Howard mapping at $\mathfrak{Q}^{(2)}\circ^{1\mathbf{Pth}_{\boldsymbol{\mathcal{A}}^{(2)}}}_{s}\mathfrak{P}^{(2)}$.

The case $i=\bb{\mathfrak{P}^{(2)}}-1$ follows.

If~(2.2.3), i.e., if $\mathfrak{Q}^{(2)}\circ^{1\mathbf{Pth}_{\boldsymbol{\mathcal{A}}^{(2)}}}_{s}
\mathfrak{P}^{(2)}$ is a head-constant echelonless second-order path that is not coherent and $i\in[ \bb{\mathfrak{P}^{(2)}}, 
\bb{\mathfrak{Q}^{(2)}\circ^{1\mathbf{Pth}_{\boldsymbol{\mathcal{A}}^{(2)}}}_{s}\mathfrak{P}^{(2)}}-1
]$ is the greatest index for which $(\mathfrak{Q}^{(2)}\circ^{1\mathbf{Pth}_{\boldsymbol{\mathcal{A}}^{(2)}}}_{s}\mathfrak{P}^{(2)})^{0,i}$ is coherent, then, regarding the second-order paths $\mathfrak{P}^{(2)}$ and $\mathfrak{Q}^{(2)}$ we have that
\begin{itemize}
\item[(i)] $\mathfrak{P}^{(2)}$ is a  coherent head-constant echelonless second-order path;
\item[(ii)] $\mathfrak{Q}^{(2),0,i-\bb{\mathfrak{P}^{(2)}}}$ is a  coherent head-constant echelonless second-order path. Moreover, $i-\bb{\mathfrak{P}^{(2)}}\in\bb{\mathfrak{Q}^{(2)}}-1$ is the greatest index for which $\mathfrak{Q}^{(2),0,i-\bb{\mathfrak{P}^{(2)}}}$ is coherent;
\item[(iii)] $\mathfrak{Q}^{(2),i-\bb{\mathfrak{P}^{(2)}}+1,\bb{\mathfrak{Q}^{(2)}}-1}$ is a head-constant echelonless second-order path.
\end{itemize}

This case follows by a similar argument to that presented in Case (2.2.1). We leave the details for the interested reader.

This completes Case~(2.2).

If~(2.3), i.e.,  if $\mathfrak{Q}^{(2)}\circ^{1\mathbf{Pth}_{\boldsymbol{\mathcal{A}}^{(2)}}}_{s}
\mathfrak{P}^{(2)}$ is a coherent head-constant echelonless second-order path then, regarding the second-order paths $\mathfrak{Q}^{(2)}$ and $\mathfrak{P}^{(2)}$, we have that 
\begin{itemize}
\item[(i)] $\mathfrak{P}^{(2)}$ is a coherent head-constant echelonless second-order path;
\item[(ii)] $\mathfrak{Q}^{(2)}$ is a coherent head-constant echelonless second-order path.
\end{itemize}

 Therefore, for a unique word $\mathbf{s}\in S^{\star}-\{\lambda\}$ and a unique operation symbol $\tau\in\Sigma^{\boldsymbol{\mathcal{A}}}_{\mathbf{s},s}$, the family of first-order translations occurring in $\mathfrak{Q}^{(2)}\circ^{1\mathbf{Pth}_{\boldsymbol{\mathcal{A}}^{(2)}}}\mathfrak{P}^{(2)}$ is a family of first-order translations of type $\tau$.

Since $\mathfrak{Q}^{(2)}\circ^{1\mathbf{Pth}_{\boldsymbol{\mathcal{A}}^{(2)}}}\mathfrak{P}^{(2)}$ is associated to the operation symbol $\tau$, the second-order paths $\mathfrak{Q}^{(2)}$ and $\mathfrak{P}^{(2)}$ are also associated to this same operation symbol.

Let $((\mathfrak{Q}^{(2)}\circ^{1\mathbf{Pth}_{\boldsymbol{\mathcal{A}}^{(2)}}}\mathfrak{P}^{(2)})_{j})_{j\in\bb{\mathbf{s}}}$ be the family of second-order paths we can extract from $\mathfrak{Q}^{(2)}\circ^{1\mathbf{Pth}_{\boldsymbol{\mathcal{A}}^{(2)}}}\mathfrak{P}^{(2)}$ in virtue of Lemma~\ref{LDPthExtract}. Then, the value of the second-order Curry-Howard mapping at $\mathfrak{Q}^{(2)}\circ^{1\mathbf{Pth}_{\boldsymbol{\mathcal{A}}^{(2)}}}\mathfrak{P}^{(2)}$ is given by
\begin{multline*}
\mathrm{CH}^{(2)}_{s}\left(
\mathfrak{Q}^{(2)}\circ^{1\mathbf{Pth}_{\boldsymbol{\mathcal{A}}^{(2)}}}\mathfrak{P}^{(2)}
\right)
\\=
\tau^{\mathbf{T}_{\Sigma^{\boldsymbol{\mathcal{A}}^{(2)}}}(X)}
\left(\left(\mathrm{CH}^{(2)}_{s_{j}}\left(
\left(\mathfrak{Q}^{(2)}\circ^{1\mathbf{Pth}_{\boldsymbol{\mathcal{A}}^{(2)}}}\mathfrak{P}^{(2)}
\right)_{j}
\right)\right)_{j\in\bb{\mathbf{s}}}
\right).
\end{multline*}

The following equality is a direct consequence of the equality proved in the first item of this lemma
\begin{multline*}
\mathrm{CH}^{(2)}_{s}
\left(
\mathrm{ip}^{(2,X)@}_{s}
\left(
\mathrm{CH}^{(2)}_{s}\left(
\mathfrak{Q}^{(2)}
\right)
\circ^{1\mathbf{T}_{\Sigma^{\boldsymbol{\mathcal{A}}^{(2)}}}(X)}_{s}
\mathrm{CH}^{(2)}_{s}\left(
\mathfrak{P}^{(2)}
\right)
\right)
\right)
\\=
\mathrm{CH}^{(2)}_{s}
\left(
\mathrm{ip}^{(2,X)@}_{s}
\left(
\mathrm{CH}^{(2)}_{s}\left(
\mathfrak{Q}^{(2)}
\right)
\right)
\circ^{1\mathbf{Pth}_{\boldsymbol{\mathcal{A}}^{(2)}}}_{s}
\mathrm{ip}^{(2,X)@}_{s}
\left(
\mathrm{CH}^{(2)}_{s}\left(
\mathfrak{P}^{(2)}
\right)
\right)
\right)
.
\end{multline*}

We will study the nature of the second-order path $\mathfrak{R}^{(2)}$ given by
$$
\mathfrak{R}^{(2)}=
\mathrm{ip}^{(2,X)@}_{s}
\left(
\mathrm{CH}^{(2)}_{s}\left(
\mathfrak{Q}^{(2)}
\right)
\right)
\circ^{1\mathbf{Pth}_{\boldsymbol{\mathcal{A}}^{(2)}}}_{s}
\mathrm{ip}^{(2,X)@}_{s}
\left(
\mathrm{CH}^{(2)}_{s}\left(
\mathfrak{P}^{(2)}
\right)
\right).
$$

Note that
\begin{itemize}
\item[(i)]  $\mathfrak{P}^{(2)}$ is a  coherent head-constant echelonless second-order path associated to the operation symbol $\tau$.
\end{itemize}

Let $(\mathfrak{P}^{(2)}_{j})_{j\in\bb{\mathbf{s}}}$ be the family of second-order paths we can extract from $\mathfrak{P}^{(2)}$ in virtue of Lemma~\ref{LDPthExtract}.  Then, according to Definition~\ref{DDCH}, the value of the second-order Curry-Howard mapping at $\mathfrak{P}^{(2)}$ is given by
$$
\mathrm{CH}^{(2)}_{s}\left(
\mathfrak{P}^{(2)}
\right)
=
\tau^{\mathbf{T}_{\Sigma^{\boldsymbol{\mathcal{A}}^{(2)}}}(X)}
\left(\left(\mathrm{CH}^{(2)}_{s_{j}}\left(
\mathfrak{P}^{(2)}_{j}
\right)\right)_{j\in\bb{\mathbf{s}}}\right).
$$

By Lemma~\ref{LDCHNEchHdC}, $\mathrm{CH}^{(2)}_{s}(\mathfrak{P}^{(2)})\in\mathcal{T}(\tau,\mathrm{T}_{\Sigma^{\boldsymbol{\mathcal{A}}^{(2)}}}(X))_{1}$, which is a subset of $\mathrm{T}_{\Sigma^{\boldsymbol{\mathcal{A}}^{(2)}}}(X)^{\mathsf{HdC}\And\mathsf{C}}_{s}$. By Proposition~\ref{PDIpDCH} $\mathrm{ip}^{(2,X)@}_{s}(\mathrm{CH}^{(2)}(\mathfrak{P}^{(2)}))$ is a second-order path in $[\mathfrak{P}^{(2)}]^{}_{s}$ then, by Lemma~\ref{LDCHNEchHdC}, we have that
\begin{itemize}
\item[(i)] $\mathrm{ip}^{(2,X)@}_{s}(\mathrm{CH}^{(2)}(\mathfrak{P}^{(2)}))$ is a coherent head-constant  echelonless second-order path associated to the operation symbol $\tau$.
\end{itemize}

Let $(\mathfrak{P}'^{(2)}_{j})_{j\in\bb{\mathbf{s}}}$ be the family of second-order paths we can extract from $\mathrm{ip}^{(2,X)@}_{s}(\mathrm{CH}^{(2)}(\mathfrak{P}^{(2)}))$ in virtue of Lemma~\ref{LDPthExtract}. Then, according to Definition~\ref{DDCH}, the value of the second-order Curry-Howard mapping at $\mathrm{ip}^{(2,X)@}_{s}(\mathrm{CH}^{(2)}(\mathfrak{P}^{(2)}))$ is given by
$$
\mathrm{CH}^{(2)}_{s}\left(
\mathrm{ip}^{(2,X)@}_{s}\left(
\mathrm{CH}^{(2)}\left(
\mathfrak{P}^{(2)}
\right)\right)\right)
=
\tau^{\mathbf{T}_{\Sigma^{\boldsymbol{\mathcal{A}}^{(2)}}}(X)}
\left(\left(\mathrm{CH}^{(2)}_{s_{j}}\left(
\mathfrak{P}'^{(2)}_{j}
\right)\right)_{j\in\bb{\mathbf{s}}}\right).
$$

By Proposition~\ref{PDIpDCH} $\mathrm{ip}^{(2,X)@}_{s}(\mathrm{CH}^{(2)}(\mathfrak{P}^{(2)}))$ is a second-order path in $[\mathfrak{P}^{(2)}]^{}_{s}$. Therefore, we have that, for every $j\in\bb{\mathbf{s}}$, it is the case that 
$$
\left(\mathfrak{P}^{(2)}_{j},
\mathfrak{P}'^{(2)}_{j}
\right)\in\mathrm{Ker}\left(\mathrm{CH}^{(2)}\right)_{s_{j}}.
$$

On the other hand,
\begin{itemize}
\item[(ii)]  $\mathfrak{Q}^{(2)}$ is a coherent head-constant  echelonless second-order path associated to the operation symbol $\tau$.
\end{itemize}

Let $(\mathfrak{Q}^{(2)}_{j})_{j\in\bb{\mathbf{s}}}$ be the family of second-order paths we can extract from $\mathfrak{Q}^{(2)}$ in virtue of Lemma~\ref{LDPthExtract}.  Then, according to Definition~\ref{DDCH}, the value of the second-order Curry-Howard mapping at $\mathfrak{Q}^{(2)}$ is given by
$$
\mathrm{CH}^{(2)}_{s}\left(
\mathfrak{Q}^{(2)}
\right)
=
\tau^{\mathbf{T}_{\Sigma^{\boldsymbol{\mathcal{A}}^{(2)}}}(X)}
\left(\left(\mathrm{CH}^{(2)}_{s_{j}}\left(
\mathfrak{Q}^{(2)}_{j}
\right)\right)_{j\in\bb{\mathbf{s}}}\right).
$$

By Lemma~\ref{LDCHNEchHdC}, $\mathrm{CH}^{(2)}_{s}(\mathfrak{Q}^{(2)})\in\mathcal{T}(\tau,\mathrm{T}_{\Sigma^{\boldsymbol{\mathcal{A}}^{(2)}}}(X))_{1}$, which is a subset of $\mathrm{T}_{\Sigma^{\boldsymbol{\mathcal{A}}}}(X)^{\mathsf{HdC}\And\mathsf{C}}_{s}$. By Proposition~\ref{PDIpDCH}, $\mathrm{ip}^{(2,X)@}_{s}(\mathrm{CH}^{(2)}(\mathfrak{Q}^{(2)}))$ is a second-order path in $[\mathfrak{Q}^{(2)}]^{}_{s}$ then, by Lemma~\ref{LDCHNEchHdC}, we have that
\begin{itemize}
\item[(ii)] $\mathrm{ip}^{(2,X)@}_{s}(\mathrm{CH}^{(2)}(\mathfrak{Q}^{(2)}))$ is a coherent head-constant  echelonless second-order path associated to the operation symbol $\tau$.
\end{itemize}

Let $(\mathfrak{Q}'^{(2)}_{j})_{j\in\bb{\mathbf{s}}}$ be the family of second-order paths we can extract from $\mathrm{ip}^{(2,X)@}_{s}(\mathrm{CH}^{(2)}(\mathfrak{Q}^{(2)}))$ in virtue of Lemma~\ref{LDPthExtract}. Then, according to Definition~\ref{DDCH}, the value of the second-order Curry-Howard mapping at $\mathrm{ip}^{(2,X)@}_{s}(\mathrm{CH}^{(2)}(\mathfrak{Q}^{(2)}))$ is given by
$$
\mathrm{CH}^{(2)}_{s}\left(
\mathrm{ip}^{(2,X)@}_{s}(\mathrm{CH}^{(2)}(\mathfrak{Q}^{(2)}))
\right)
=
\tau^{\mathbf{T}_{\Sigma^{\boldsymbol{\mathcal{A}}^{(2)}}}(X)}
\left(\left(\mathrm{CH}^{(2)}_{s_{j}}\left(
\mathfrak{Q}'^{(2)}_{j}
\right)\right)_{j\in\bb{\mathbf{s}}}\right).
$$

By Proposition~\ref{PDIpDCH}, $\mathrm{ip}^{(2,X)@}_{s}(\mathrm{CH}^{(2)}(\mathfrak{Q}^{(2)}))$  is a second-order path in $[\mathfrak{Q}^{(2)}]^{}_{s}$. Therefore, we have that, for every $j\in\bb{\mathbf{s}}$, it is the case that 
$$
\left(\mathfrak{Q}^{(2)}_{j},
\mathfrak{Q}'^{(2)}_{j}
\right)\in\mathrm{Ker}\left(\mathrm{CH}^{(2)}\right)_{s_{j}}.
$$

Taking into account items (i) and (ii) and the fact that $\mathfrak{Q}^{(2)}\circ^{1\mathbf{Pth}_{\boldsymbol{\mathcal{A}}^{(2)}}}_{s}
\mathfrak{P}^{(2)}$ is a coherent head-constant echelonless second-order path, we conclude in virtue of Proposition~\ref{PDIpDCH} and Lemma~\ref{LTech} that the 
the second-order path $\mathfrak{R}^{(2)}$ given by
$$
\mathfrak{R}^{(2)}=
\mathrm{ip}^{(2,X)@}_{s}\left(
\mathrm{CH}^{(2)}_{s}\left(
\mathfrak{Q}^{(2)}
\right)\right)
\circ^{1\mathbf{Pth}_{\boldsymbol{\mathcal{A}}^{(2)}}}_{s}
\mathrm{ip}^{(2,X)@}_{s}\left(
\mathrm{CH}^{(2)}_{s}\left(
\mathfrak{P}^{(2)}
\right)\right)
$$
is a  coherent head-constant echelonless second-order path  associated to the operation symbol $\tau$.

Let $(\mathfrak{R}^{(2)}_{j})_{j\in\bb{\mathbf{s}}}$ be the family of second-order paths we can extract from $\mathfrak{R}^{(2)}$ in virtue of Lemma~\ref{LDPthExtract}. For every $j\in\bb{\mathbf{s}}$, the following chain of equalities holds
\allowdisplaybreaks
\begin{align*}
\mathfrak{R}^{(2)}_{j}
&=
\mathfrak{Q}'^{(2)}_{j}
\circ^{1\mathbf{Pth}_{\boldsymbol{\mathcal{A}}^{(2)}}}_{s_{j}}
\mathfrak{P}'^{(2)}_{j}
\tag{1}
\\&=
\mathrm{ip}^{(2,X)@}_{s_{j}}\left(
\mathrm{CH}^{(2)}_{s_{j}}\left(
\mathfrak{Q}^{(2)}_{j}
\right)\right)
\circ^{1\mathbf{Pth}_{\boldsymbol{\mathcal{A}}^{(2)}}}_{s_{j}}
\mathrm{ip}^{(2,X)@}_{s_{j}}\left(
\mathrm{CH}^{(2)}_{s_{j}}\left(
\mathfrak{P}^{(2)}_{j}
\right)\right).
\tag{2}
\end{align*}

The first equality follows from the proof of Lemma~\ref{LDPthExtract} and the description of $\mathfrak{R}^{(2)}$; finally, the last equality follows from the description of $\mathfrak{R}^{(2)}$ as a $1$-composition of normalized paths. Thus, the $j$-th component we can extract from the normalized second-order path $\mathrm{ip}^{(2,X)@}_{s}(\mathrm{CH}^{(2)}_{s}(\mathfrak{Q}^{(2)}))$ is the normalization of the $j$-th component we can extract from $\mathfrak{Q}^{(2)}$. Analogously for $\mathrm{ip}^{(2,X)@}_{s}(\mathrm{CH}^{(2)}_{s}(\mathfrak{P}^{(2)}))$.

Also from the proof of Lemma~\ref{LDPthExtract}, we also have that, for every $j\in\bb{\mathbf{s}}$,
$$
\left(\mathfrak{Q}^{(2)}
\circ^{1\mathbf{Pth}_{\boldsymbol{\mathcal{A}}^{(2)}}}_{s}
\mathfrak{P}^{(2)}\right)_{j}
=
\mathfrak{Q}^{(2)}_{j}
\circ^{1\mathbf{Pth}_{\boldsymbol{\mathcal{A}}^{(2)}}}_{s_{j}}
\mathfrak{P}^{(2)}_{j}.
$$

Therefore, we are in position to conclude the initial discussion at the beginning of this subcase. The following chain of equalities holds
\begin{flushleft}
$\mathrm{CH}^{(2)}_{s}
\left(
\mathrm{ip}^{(2,X)@}_{s}
\left(
\mathrm{CH}^{(2)}_{s}\left(
\mathfrak{Q}^{(2)}
\right)
\circ^{1\mathbf{T}_{\Sigma^{\boldsymbol{\mathcal{A}}^{(2)}}}(X)}_{s}
\mathrm{CH}^{(2)}_{s}\left(
\mathfrak{P}^{(2)}
\right)
\right)
\right)
$
\allowdisplaybreaks
\begin{align*}
\quad&=
\mathrm{CH}^{(2)}_{s}
\left(
\mathrm{ip}^{(2,X)@}_{s}
\left(
\mathrm{CH}^{(2)}_{s}\left(
\mathfrak{Q}^{(2)}
\right)
\right)
\circ^{1\mathbf{Pth}_{\boldsymbol{\mathcal{A}}^{(2)}}}_{s}
\mathrm{ip}^{(2,X)@}_{s}
\left(
\mathrm{CH}^{(2)}_{s}\left(
\mathfrak{P}^{(2)}
\right)
\right)
\right)
\tag{1}
\\&=
\tau^{\mathbf{T}_{\Sigma^{\boldsymbol{\mathcal{A}}^{(2)}}}(X)}
\Bigg(\left(\mathrm{CH}^{(2)}_{s_{j}}\left(
\mathrm{ip}^{(2,X)@}_{s_{j}}\left(
\mathrm{CH}^{(2)}_{s_{j}}\left(
\mathfrak{Q}^{(2)}_{j}
\right)\right)
\circ^{1\mathbf{Pth}_{\boldsymbol{\mathcal{A}}^{(2)}}}_{s_{j}}
\right.\right.
\\&\qquad\qquad\qquad\qquad\qquad\qquad\qquad\qquad
\left.\left.
\mathrm{ip}^{(2,X)@}_{s_{j}}\left(
\mathrm{CH}^{(2)}_{s_{j}}\left(
\mathfrak{P}^{(2)}_{j}
\right)\right)
\right)
\right)_{j\in\bb{\mathbf{s}}}
\Bigg)
\tag{2}
\\&=
\tau^{\mathbf{T}_{\Sigma^{\boldsymbol{\mathcal{A}}^{(2)}}}(X)}
\Bigg(\left(\mathrm{CH}^{(2)}_{s_{j}}\left(
\mathrm{ip}^{(2,X)@}_{s_{j}}\left(
\mathrm{CH}^{(2)}_{s_{j}}\left(
\mathfrak{Q}^{(2)}_{j}
\right)
\circ^{1\mathbf{T}_{\Sigma^{\boldsymbol{\mathcal{A}}^{(2)}}}(X)}_{s_{j}}
\right.\right.\right.
\\&\qquad\qquad\qquad\qquad\qquad\qquad\qquad\qquad\qquad\qquad\quad\,\,\,\,
\mathrm{CH}^{(2)}_{s_{j}}\left(
\mathfrak{P}^{(2)}_{j}
\right)\bigg)\bigg)
\bigg)_{j\in\bb{\mathbf{s}}}
\Bigg)
\tag{3}
\\&=
\tau^{\mathbf{T}_{\Sigma^{\boldsymbol{\mathcal{A}}^{(2)}}}(X)}
\left(\left(\mathrm{CH}^{(2)}\left(
\mathfrak{Q}^{(2)}_{j}
\circ^{1\mathbf{Pth}_{\boldsymbol{\mathcal{A}}^{(2)}}}_{s_{j}}
\mathfrak{P}^{(2)}_{j}
\right)
\right)_{j\in\bb{\mathbf{s}}}
\right)
\tag{4}
\\&=
\mathrm{CH}^{(2)}_{s}\left(
\mathfrak{Q}^{(2)}
\circ^{1\mathbf{Pth}_{\boldsymbol{\mathcal{A}}^{(2)}}}_{s}
\mathfrak{P}^{(2)}
\right).
\tag{5}
\end{align*}
\end{flushleft}

The first equality was already proven at the beginning of this case; the second equality simply unravels the second-order Curry-Howard mapping at the second-order path $\mathfrak{R}^{(2)}$; the third equality follows from the fact that, by Definition~\ref{DDIp}, $\mathrm{ip}^{(2,X)@}$ is a $\Sigma^{\boldsymbol{\mathcal{A}}^{(2)}}$-homomorphism; the fourth equality follows from induction. Let us note that the pair $(\mathfrak{Q}^{(2)}_{j}\circ^{1\mathbf{Pth}_{\boldsymbol{\mathcal{A}}^{(2)}}}_{s_{j}}\mathfrak{P}^{(2)}_{j},s_{j})$ $\prec_{\mathbf{Pth}_{\boldsymbol{\mathcal{A}}^{(2)}}}$-precedes $(\mathfrak{Q}^{(2)}\circ^{1\mathbf{Pth}_{\boldsymbol{\mathcal{A}}^{(2)}}}_{s}\mathfrak{P}^{(2)},s)$, thus we have that
\begin{multline*}
\mathrm{CH}^{(2)}_{s_{j}}
\left(
\mathrm{ip}^{(2,X)@}_{s_{j}}
\left(
\mathrm{CH}^{(2)}_{s_{j}}\left(
\mathfrak{Q}^{(2)}_{j}
\right)
\circ^{1\mathbf{T}_{\Sigma^{\boldsymbol{\mathcal{A}}^{(2)}}}(X)}_{s_{j}}
\mathrm{CH}^{(2)}_{s_{j}}\left(
\mathfrak{P}^{(2)}_{j}
\right)
\right)
\right)
\\ =
\mathrm{CH}^{(2)}_{s_{j}}\left(
\mathfrak{Q}^{(2)}_{j}
\circ^{1\mathbf{Pth}_{\boldsymbol{\mathcal{A}}^{(2)}}}_{s_{j}}
\mathfrak{P}^{(2)}_{j}
\right);
\end{multline*}
finally, the last equality recovers the value of the second-order Curry-Howard mapping at $\mathfrak{Q}^{(2)}\circ^{1\mathbf{Pth}_{\boldsymbol{\mathcal{A}}^{(2)}}}_{s}\mathfrak{P}^{(2)}$.

Case (2.3) follows.

This completes the proof.
\end{proof}

\section{
\texorpdfstring
{The congruence $\Theta^{[2]}$ on $\mathbf{T}_{\Sigma^{\boldsymbol{\mathcal{A}}^{(2)}}}(X)$}
{A  congruence for the second-order free algebra}
}

We next show that the properties that hold for the relation $\Theta^{(2)}$, set up in Definition~\ref{DDTheta} also hold for the smallest $\Sigma^{\boldsymbol{\mathcal{A}}^{(2)}}$-congruence generated by it. In this regard, we recover the notation introduced in Definition~\ref{DCongOpInt}.

\begin{restatable}{definition}{DDThetaCong}
\label{DDThetaCong}
\index{Theta!second-order!$\Theta^{[2]}$}
 We denote by $\Theta^{[2]}$ the smallest $\Sigma^{\boldsymbol{\mathcal{A}}^{(2)}}$-congruence on  
$\mathbf{T}_{\Sigma^{\boldsymbol{\mathcal{A}}^{(2)}}}(X)$ containing $\Theta^{(2)}$, i.e.,
$$\Theta^{[2]}
=\mathrm{Cg}_{\mathbf{T}_{\Sigma^{\boldsymbol{\mathcal{A}}^{(2)}}}(X)}
\left(\Theta^{(2)}
\right).
$$
\end{restatable}

\begin{remark}\label{RDThetaCong} Let $\mathrm{C}_{\Sigma^{\boldsymbol{\mathcal{A}}^{(2)}}}$ stand for the operator $\mathrm{C}_{\mathbf{T}_{\Sigma^{\boldsymbol{\mathcal{A}}^{(2)}}}(X)}$ on $\mathrm{T}_{\Sigma^{\boldsymbol{\mathcal{A}}^{(2)}}}(X)\times \mathrm{T}_{\Sigma^{\boldsymbol{\mathcal{A}}^{(2)}}}(X)$ (see Definition~\ref{DCongOpC} for the general case). Let us also recall, from Definition~\ref{DCongOpC}, that if $\Phi\subseteq \mathrm{T}_{\Sigma^{\boldsymbol{\mathcal{A}}^{(2)}}}(X)\times \mathrm{T}_{\Sigma^{\boldsymbol{\mathcal{A}}^{(2)}}}(X)$, then 
$$
\mathrm{C}_{\Sigma^{\boldsymbol{\mathcal{A}}^{(2)}}}(\Phi)
=
(\Phi\circ\Phi)
\cup
\left(
\bigcup_{\gamma\in\Sigma^{\boldsymbol{\mathcal{A}}^{(2)}}_{\neq\lambda,s}}
\gamma^{\mathbf{T}_{\Sigma^{\boldsymbol{\mathcal{A}}^{(2)}}}(X)}
\times
\gamma^{\mathbf{T}_{\Sigma^{\boldsymbol{\mathcal{A}}^{(2)}}}(X)}
\left[
\Phi_{\mathrm{ar}(\gamma)}
\right]
\right)_{s\in S}.
$$
Moreover, for the family $(\mathrm{C}^{n}_{\Sigma^{\boldsymbol{\mathcal{A}}^{(2)}}}(\Theta^{(2)}))_{n\in\mathbb{N}}$ in $\mathrm{Sub}(\mathrm{T}_{\Sigma^{\boldsymbol{\mathcal{A}}^{(2)}}}(X)^{2}
)
$, defined recursively as follows:
\allowdisplaybreaks
\begin{align*}
\mathrm{C}^{0}_{\Sigma^{\boldsymbol{\mathcal{A}}^{(2)}}}
\left(
\Theta^{(2)}
\right)
&=\left(\Theta^{(2)}\right)
\cup
\left(
\Theta^{(2)}
\right)^{-1}
\cup\Delta_{\mathrm{T}_{\Sigma^{\boldsymbol{\mathcal{A}}^{(2)}}}(X)},
\\
\mathrm{C}^{n+1}_{\Sigma^{\boldsymbol{\mathcal{A}}^{(2)}}}\left(
\Theta^{(2)}
\right)
&=\mathrm{C}_{\Sigma^{\boldsymbol{\mathcal{A}}^{(2)}}}\left(
\mathrm{C}^{n+1}_{\Sigma^{\boldsymbol{\mathcal{A}}^{(2)}}}\left(
\Theta^{(2)}
\right)\right),\,n\geq 0,
\end{align*}
we have, by Proposition~\ref{PCongOpC}, that 
$$\mathrm{C}^{\omega}_{\Sigma^{\boldsymbol{\mathcal{A}}^{(2)}}}\left(
\Theta^{(2)}
\right)
=
\bigcup_{n\in\mathbb{N}}\mathrm{C}^{n}_{\Sigma^{\boldsymbol{\mathcal{A}}^{(2)}}}\left(
\Theta^{(2)}\right)=\Theta^{[2]},$$ is equal to the smallest $\Sigma^{\boldsymbol{\mathcal{A}}^{(2)}}$-congruence on $\mathbf{T}_{\Sigma^{\boldsymbol{\mathcal{A}}^{(2)}}}(X)$ containing $\Theta^{(2)}$.
\end{remark}

\begin{restatable}{proposition}{PDThetaCongDCatAlg}
\label{PDThetaCongDCatAlg}
\index{terms!second-order!$\mathrm{T}_{\Sigma^{\boldsymbol{\mathcal{A}}^{(2)}}}(X)/{\Theta^{[2]}}$}
\index{terms!second-order!$[P]_{\Theta^{[2]}_{s}}$}
 The $S$-sorted set $\mathrm{T}_{\Sigma^{\boldsymbol{\mathcal{A}}^{(2)}}}(X)/{\Theta^{[2]}}$ is equipped, in a natural way, with a structure of many-sorted $\Sigma^{\boldsymbol{\mathcal{A}}^{(2)}}$-algebra. We denote by $\mathbf{T}_{\Sigma^{\boldsymbol{\mathcal{A}}^{(2)}}}(X)/{\Theta^{[2]}}$ the corresponding $\Sigma^{\boldsymbol{\mathcal{A}}^{(2)}}$-algebra.

\index{projection!second-order!$\mathrm{pr}^{\Theta^{[2]}}$}
We define the projection $\mathrm{pr}^{\Theta^{[2]}}$ from $\mathrm{T}_{\Sigma^{\boldsymbol{\mathcal{A}}^{(2)}}}(X)$ to $\mathrm{T}_{\Sigma^{\boldsymbol{\mathcal{A}}^{(2)}}}(X)/{\Theta^{[2]}}$ to be the many-sorted mapping that, for every sort $s\in S$, maps a term $P$ in $\mathrm{T}_{\Sigma^{\boldsymbol{\mathcal{A}}^{(2)}}}(X)_{s}$ to $[P]_{\Theta^{[2]}_{s}}$, its equivalence class under the $\Sigma^{\boldsymbol{\mathcal{A}}^{(2)}}$-congruence $\Theta^{[2]}$.

Note that $\mathrm{pr}^{\Theta^{[2]}}$, that is,
$$
\mathrm{pr}^{\Theta^{[2]}}\colon
\mathbf{T}_{\Sigma^{\boldsymbol{\mathcal{A}}^{(2)}}}(X)
\mor
\mathbf{T}_{\Sigma^{\boldsymbol{\mathcal{A}}^{(2)}}}(X)/{\Theta^{[2]}}
$$
is a surjective $\Sigma^{\boldsymbol{\mathcal{A}}^{(2)}}$-homomorphism from $\mathbf{T}_{\Sigma^{\boldsymbol{\mathcal{A}}^{(2)}}}(X)$ to $\mathbf{T}_{\Sigma^{\boldsymbol{\mathcal{A}}^{(2)}}}(X)/{\Theta^{[2]}}$.
\end{restatable}

We next introduce the reductions of the $\Sigma^{\boldsymbol{\mathcal{A}}^{(2)}}$-algebra $\mathbf{T}_{\Sigma^{\boldsymbol{\mathcal{A}}^{(2)}}}(X)/{\Theta^{[2]}}$ to layers $0$ and $1$.

\begin{definition} We will denote by $\mathbf{T}^{(1,2)}_{\Sigma^{\boldsymbol{\mathcal{A}}^{(2)}}}(X)/{\Theta^{[2]}}$ the $\Sigma^{\boldsymbol{\mathcal{A}}}$-algebra 
$$\left(\mathbf{in}^{\Sigma,(1,2)}
\right)\left(\mathbf{T}_{\Sigma^{\boldsymbol{\mathcal{A}}^{(2)}}}(X)/{\Theta^{[2]}}\right).$$

We will call $\mathbf{T}^{(1,2)}_{\Sigma^{\boldsymbol{\mathcal{A}}^{(2)}}}(X)/{\Theta^{[2]}}$ the \emph{$\Sigma^{\boldsymbol{\mathcal{A}}}$-reduct} of the many-sorted $\Sigma^{\boldsymbol{\mathcal{A}}^{(2)}}$-algebra $\mathbf{T}_{\Sigma^{\boldsymbol{\mathcal{A}}^{(2)}}}(X)/{\Theta^{[2]}}$. In this regard see Definition~\ref{DDURed} and Remark~\ref{RDURed}.

We will denote by $\mathbf{T}^{(0,2)}_{\Sigma^{\boldsymbol{\mathcal{A}}^{(2)}}}(X)/{\Theta^{[2]}}$ the $\Sigma$-algebra 
$$\left(\mathbf{in}^{\Sigma,(0,2)}
\right)\left(\mathbf{T}_{\Sigma^{\boldsymbol{\mathcal{A}}^{(2)}}}(X)/{\Theta^{[2]}}\right).$$

We will call $\mathbf{T}^{(0,2)}_{\Sigma^{\boldsymbol{\mathcal{A}}^{(2)}}}(X)/{\Theta^{[2]}}$ the \emph{$\Sigma$-reduct} of the many-sorted $\Sigma^{\boldsymbol{\mathcal{A}}^{(2)}}$-algebra $\mathbf{T}_{\Sigma^{\boldsymbol{\mathcal{A}}^{(2)}}}(X)/{\Theta^{[2]}}$. In this regard see Definition~\ref{DDZRed} and Remark~\ref{RDZRed}.
\end{definition}

In the following lemmas we investigate the $\Sigma^{\boldsymbol{\mathcal{A}}^{(2)}}$-congruence $\Theta^{[2]}$ on $\mathbf{T}_{\Sigma^{\boldsymbol{\mathcal{A}}^{(2)}}}(X)$. In particular, we will show how this congruence can be used to describe those terms that when mapped to $\mathrm{F}_{\Sigma^{\boldsymbol{\mathcal{A}}^{(2)}}}(\mathbf{Pth}_{\boldsymbol{\mathcal{A}}^{(2)}})$ under the action of $\mathrm{ip}^{(2,X)@}$ retrieve a second-order path.

We start with a lemma that will be of great importance later on. Actually, it states that if a term is such that its image under $\mathrm{ip}^{(2,X)@}$ retrieves a second-order path, then all its subterms have the same property.

\begin{restatable}{lemma}{LDThetaCongSub}
\label{LDThetaCongSub} Let $s$ be a sort in $S$ and $P$ a term in $\mathrm{T}_{\Sigma^{\boldsymbol{\mathcal{A}}^{(2)}}}(X)_{s}$. If $\mathrm{ip}^{(2,X)@}_{s}(P)$ is a second-order path in $\mathrm{Pth}_{\boldsymbol{\mathcal{A}}^{(2)},s}$, then $\mathrm{ip}^{(2,X)@}[\mathrm{Subt}_{\Sigma^{\boldsymbol{\mathcal{A}}^{(2)}}}(P)]$ is a subset of $\mathrm{Pth}_{\boldsymbol{\mathcal{A}}^{(2)}}$.
\end{restatable}
\begin{proof}
We prove it by induction on the height of $P$.

\textsf{Base step of the induction.}

If $P$ is a term of height $0$, that is, if $P\in\mathrm{B}^{0}_{\Sigma^{\boldsymbol{\mathcal{A}}^{(2)}}}(X)_{s}$, then
$$
\mathrm{Subt}_{\Sigma^{\boldsymbol{\mathcal{A}}^{(2)}}}(P)=\delta^{P,s}.
$$
That is, the unique subterm of $P$ is $P$ itself. The statement trivially holds.

\textsf{Inductive step of the induction.}

Assume that the statement holds for terms of height up to $n$, that is, for every sort $s\in S$, if $P\in\mathrm{B}^{n}_{\Sigma^{\boldsymbol{\mathcal{A}}^{(2)}}}(X)_{s}$ is such that $\mathrm{ip}^{(2,X)@}_{s}(P)$ is a second-order path then $\mathrm{ip}^{(2,X)@}[\mathrm{Subt}_{\Sigma^{\boldsymbol{\mathcal{A}}^{(2)}}}(P)]$ is a subset of second-order paths in $\mathrm{Pth}_{\boldsymbol{\mathcal{A}}^{(2)}}$.

Now, let $P\in\mathrm{B}^{n+1}_{\Sigma^{\boldsymbol{\mathcal{A}}^{(2)}}}(X)_{s}$ be a term of height $n+1$. Then there exists a unique word $\mathbf{s}\in S^{\star}-\{\lambda\}$, a unique operation symbol $\gamma\in\Sigma^{\boldsymbol{\mathcal{A}}^{(2)}}_{\mathbf{s},s}$, and a unique family of terms $(P_{j})_{j\in\bb{\mathbf{s}}}\in\mathrm{B}^{n}_{\Sigma^{\boldsymbol{\mathcal{A}}^{(2)}}}(X)_{\mathbf{s}}$ for which $P=\gamma^{\mathbf{T}_{\Sigma^{\boldsymbol{\mathcal{A}}^{(2)}}}(X)}((P_{j})_{j\in\bb{\mathbf{s}}})$. In this case, 
$$
\mathrm{Subt}_{\Sigma^{\boldsymbol{\mathcal{A}}^{(2)}}}(P)
=
\delta^{P,s}
\cup
\bigcup_{j\in\bb{\mathbf{s}}}
\mathrm{Subt}_{\Sigma^{\boldsymbol{\mathcal{A}}^{(2)}}}
(P_{j}).
$$

Note that the following chain of equivalences holds
\begin{flushleft}
$\mathrm{ip}^{(2,X)@}_{s}(P)\in\mathrm{Pth}_{\boldsymbol{\mathcal{A}}^{(2)},s}$
\allowdisplaybreaks
\begin{align*}
\qquad
&\Leftrightarrow
\mathrm{ip}^{(2,X)@}_{s}\left(
\gamma^{\mathbf{T}_{\Sigma^{\boldsymbol{\mathcal{A}}^{(2)}}}(X)}
\left(\left(P_{j}\right)_{j\in\bb{\mathbf{s}}}\right)\right)
\in\mathrm{Pth}_{\boldsymbol{\mathcal{A}}^{(2)},s}
\tag{1}
\\&\Leftrightarrow
\gamma^{\mathbf{F}_{\Sigma^{\boldsymbol{\mathcal{A}}^{(2)}}}\left(
\mathbf{Pth}_{\boldsymbol{\mathcal{A}}^{(2)}}
\right)}\left(\left(
\mathrm{ip}^{(2,X)@}_{s_{j}}\left(
P_{j}
\right)\right)_{j\in\bb{\mathbf{s}}}\right)
\in\mathrm{Pth}_{\boldsymbol{\mathcal{A}}^{(2)},s}
\tag{2}
\\&\Leftrightarrow
\left(\mathrm{ip}^{(2,X)@}_{s_{j}}\left(
P_{j}\right)
\right)_{j\in\bb{\mathbf{s}}}
\in\mathrm{Dom}\left(
\gamma^{\mathbf{Pth}_{\boldsymbol{\mathcal{A}}^{(2)}}}
\right).
\tag{3}
\end{align*}
\end{flushleft}

The first equivalence unravels the definition of $P$; the second equation follows from the fact that, according to Definition~\ref{DDIp}, $\mathrm{ip}^{(2,X)@}$ is a $\Sigma^{\boldsymbol{\mathcal{A}}^{(2)}}$-homomorphism; the third  equivalence follows from the fact that $\gamma^{\mathbf{F}_{\Sigma^{\boldsymbol{\mathcal{A}}^{(2)}}}(
\mathbf{Pth}_{\boldsymbol{\mathcal{A}}^{(2)}}
)}((
\mathrm{ip}^{(2,X)@}_{s_{j}}(
P_{j}
))_{j\in\bb{\mathbf{s}}})$ is a second-order path in $\mathrm{Pth}_{\boldsymbol{\mathcal{A}}^{(2)},s}$ if, and only if, it is the case that $(\mathrm{ip}^{(2,X)@}_{s_{j}}(
P_{j})
)_{j\in\bb{\mathbf{s}}}$ is a family of second-order paths in the domain of $\gamma^{\mathbf{Pth}_{\boldsymbol{\mathcal{A}}^{(2)}}}$.

This completes the proof.
\end{proof}

The following proposition guarantees that the lifting of the congruence $\Theta^{[1]}$ by means of $\eta^{(2,1)\sharp}\times \eta^{(2,1)\sharp}$ is included in $\Theta^{[2]}$.

\begin{restatable}{proposition}{PDThetaCongDU}
\label{PDThetaCongDU} The following inclusion holds
\[
\eta^{(2,1)\sharp}\times\eta^{(2,1)\sharp}\left[
\Theta^{[1]}
\right]
\subseteq 
\Theta^{[2]}.
\]
\end{restatable}
\begin{proof}
Let us recall that, according to Definition~\ref{DThetaCong}, $\Theta^{[1]}$ was defined as the smallest $\Sigma^{\boldsymbol{\mathcal{A}}}$-congruence on $\mathbf{T}_{\Sigma^{\boldsymbol{\mathcal{A}}}}(X)$ containing $\Theta^{(1)}$. 

The statement follows from the fact that, according to Proposition~\ref{PDThetaDU}, the following inclusion holds 
\[\eta^{(2,1)\sharp}\times \eta^{(2,1)\sharp}\left[
\Theta^{(1)}
\right]
\subseteq \Theta^{(2)}.
\]

Let us also recall that, according to Proposition~\ref{PDUEmb}, $\eta^{(2,1)\sharp}$ is a $\Sigma^{\boldsymbol{\mathcal{A}}}$-homomorphism.

This completes the proof.
\end{proof}

The following proposition ensures that the composition of the second-order Curry-Howard mapping with the projection with respect to the $\Sigma^{\boldsymbol{\mathcal{A}}^{(2)}}$-congruence $\Theta^{[2]}$ determines a $\Sigma^{\boldsymbol{\mathcal{A}}^{(2)}}$-homomorphism from $\mathbf{Pth}_{\boldsymbol{\mathcal{A}}^{(2)}}$ to $\mathbf{T}_{\Sigma^{\boldsymbol{\mathcal{A}}^{(2)}}}(X)/{\Theta^{[2]}}$.

\begin{restatable}{proposition}{PDThetaDCH}
\label{PDThetaDCH} The mapping given by the composition
\[
\mathrm{pr}^{\Theta^{[2]}}\circ \mathrm{CH}^{(2)}
\colon\mathbf{Pth}_{\boldsymbol{\mathcal{A}}^{(2)}}
\mor
\mathbf{T}_{\Sigma^{\boldsymbol{\mathcal{A}}^{(2)}}}(X)/{\Theta^{[2]}}
\]
is a $\Sigma^{\boldsymbol{\mathcal{A}}^{(2)}}$-homomorphism.
\end{restatable} 
\begin{proof}
We prove that this composition mapping is compatible with every operation symbol in $\Sigma^{\boldsymbol{\mathcal{A}}^{(2)}}$.

\textsf{The mapping $\mathrm{pr}^{\Theta^{[2]}}\circ \mathrm{CH}^{(2)}$ is a $\Sigma$-homomorphism.}

Let us note that, according to Proposition~\ref{PDCHHom}, the second-order Curry-Howard mapping is a $\Sigma$-homomorphism from $\mathbf{Pth}_{\boldsymbol{\mathcal{A}}^{(2)}}^{(0,2)}$ to $\mathbf{T}_{\Sigma^{\boldsymbol{\mathcal{A}}^{(2)}}}^{(0,2)}(X)$. Note that $\mathrm{pr}^{\Theta^{[2]}}$ is a $\Sigma^{\boldsymbol{\mathcal{A}}^{(2)}}$-homomorphism according to Proposition~\ref{PDThetaCongDCatAlg}. Therefore, we conclude that the composition $\mathrm{pr}^{\Theta^{[2]}}\circ \mathrm{CH}^{(2)}$ is a $\Sigma$-homomorphism from $\mathbf{Pth}_{\boldsymbol{\mathcal{A}}^{(2)}}^{(0,2)}$ to $\mathbf{T}_{\Sigma^{\boldsymbol{\mathcal{A}}^{(2)}}}^{(0,2)}(X)/{\Theta^{[2]}}$.

\textsf{The mapping $\mathrm{pr}^{\Theta^{[2]}}\circ \mathrm{CH}^{(2)}$ is compatible with the rewrite rules.}

Let $s$ be a sort in $S$ and let $\mathfrak{p}$ be a rewrite rule in $\mathcal{A}_{s}$.
The following chain of equalities holds.
\allowdisplaybreaks
\begin{align*}
\mathrm{pr}^{\Theta^{[2]}}_{s}\left(
\mathrm{CH}^{(2)}_{s}\left(
\mathfrak{p}^{\mathbf{Pth}_{\boldsymbol{\mathcal{A}}^{(2)}}}
\right)
\right)&=
\mathrm{pr}^{\Theta^{[2]}}_{s}\left(
\mathfrak{p}^{\mathbf{T}_{\Sigma^{\boldsymbol{\mathcal{A}}^{(2)}}}(X)}
\right)
\tag{1}
\\&=
\mathfrak{p}^{\mathbf{T}_{\Sigma^{\boldsymbol{\mathcal{A}}^{(2)}}}(X)/{\Theta^{[2]}}}.
\tag{2}
\end{align*}

In the just stated chain of equalities, the first equality follows from Proposition~\ref{PDCHA} and the second equality follows from Proposition~\ref{PDThetaCongDCatAlg}.

\textsf{The mapping $\mathrm{pr}^{\Theta^{[2]}}\circ \mathrm{CH}^{(2)}$ is compatible with the $0$-source.}

Let $s$ be a sort in $S$ and let us consider the $0$-source operation symbol $\mathrm{sc}^{0}_{s}$ in $\Sigma^{\boldsymbol{\mathcal{A}}^{(2)}}_{s,s}$. Let $\mathfrak{P}^{(2)}$ be a second-order path in $\mathrm{Pth}_{\boldsymbol{\mathcal{A}}^{(2)},s}$. 

The following chain of equalities holds
\allowdisplaybreaks
\begin{align*}
\mathrm{pr}^{\Theta^{[2]}}_{s}\left(
\mathrm{CH}^{(2)}_{s}\left(
\mathrm{sc}^{0\mathbf{Pth}_{\boldsymbol{\mathcal{A}}^{(2)}}}_{s}
\left(
\mathfrak{P}^{(2)}
\right)\right)
\right)&=
\mathrm{pr}^{\Theta^{[2]}}_{s}\left(
\mathrm{sc}^{0\mathbf{T}_{\Sigma^{\boldsymbol{\mathcal{A}}^{(2)}}}}_{s}\left(
\mathrm{CH}^{(2)}_{s}\left(
\mathfrak{P}^{(2)}
\right)\right)
\right)
\tag{1}
\\&=
\mathrm{sc}^{0\mathbf{T}_{\Sigma^{\boldsymbol{\mathcal{A}}^{(2)}}}/{\Theta^{[2]}}}_{s}
\left(
\mathrm{pr}^{\Theta^{[2]}}_{s}\left(
\mathrm{CH}^{(2)}_{s}\left(
\mathfrak{P}^{(2)}
\right)\right)
\right).
\tag{2}
\end{align*}

In the just stated chain of equalities, the first equality follows from the fact that, by Definition~\ref{DDTheta}, the pair
\[
\left(\mathrm{CH}^{(2)}_{s}\left(
\mathrm{sc}^{0\mathbf{Pth}_{\boldsymbol{\mathcal{A}}^{(2)}}}_{s}\left(
\mathfrak{P}^{(2)}
\right)\right),
\mathrm{sc}^{0\mathbf{T}_{\Sigma^{\boldsymbol{\mathcal{A}}^{(2)}}}(X)}_{s}\left(
\mathrm{CH}^{(2)}_{s}\left(
\mathfrak{P}^{(2)}
\right)\right)
\right)
\]
is in $\Theta^{(2)}_{s}$; finally, the second equality follows from Proposition~\ref{PDThetaCongDCatAlg}.

Hence, $\mathrm{pr}^{\Theta^{[2]}}\circ \mathrm{CH}^{(2)}$ is compatible with the $0$-source operation.

\textsf{The mapping $\mathrm{pr}^{\Theta^{[2]}}\circ \mathrm{CH}^{(2)}$ is compatible with the $0$-target.}

Let $s$ be a sort in $S$ and let us consider the $0$-target operation symbol $\mathrm{tg}^{0}_{s}$ in $\Sigma^{\boldsymbol{\mathcal{A}}^{(2)}}_{s,s}$. Let $\mathfrak{P}^{(2)}$ be a second-order path in $\mathrm{Pth}_{\boldsymbol{\mathcal{A}}^{(2)},s}$, then the following equality holds
\[
\mathrm{pr}^{\Theta^{[2]}}_{s}\left(
\mathrm{CH}^{(2)}_{s}\left(
\mathrm{tg}^{0\mathbf{Pth}_{\boldsymbol{\mathcal{A}}^{(2)}}}_{s}\left(
\mathfrak{P}^{(2)}
\right)\right)\right)
=
\mathrm{tg}^{0\mathbf{T}_{\Sigma^{\boldsymbol{\mathcal{A}}^{(2)}}}(X)/{\Theta^{[2]}}}_{s}\left(
\mathrm{pr}^{\Theta^{[2]}}_{s}\left(
\mathrm{CH}^{(2)}_{s}\left(
\mathfrak{P}^{(2)}
\right)\right)\right).
\]

The proof of this case is identical to that of the $0$-source.

Hence, $\mathrm{pr}^{\Theta^{[2]}}\circ \mathrm{CH}^{(2)}$ is compatible with the $0$-target operation.

\textsf{The mapping $\mathrm{pr}^{\Theta^{[2]}}\circ \mathrm{CH}^{(2)}$ is compatible with the $0$-composition.}

Let $s$ be a sort in $S$ and let us consider the $0$-composition operation symbol $\circ^{0}_{s}$ in $\Sigma^{\boldsymbol{\mathcal{A}}^{(2)}}_{ss,s}$. Let $\mathfrak{P}^{(2)}$ and $\mathfrak{Q}^{(2)}$ be two second-order paths in $\mathrm{Pth}_{\boldsymbol{\mathcal{A}}^{(2)},s}$ satisfying that 
\[
\mathrm{sc}^{(0,2)}_{s}\left(\mathfrak{Q}^{(2)}\right)
=
\mathrm{tg}^{(0,2)}_{s}\left(\mathfrak{P}^{(2)}\right).
\]

The following chain of equalities holds
\begin{flushleft}
$\mathrm{pr}^{\Theta^{[2]}}_{s}\left(
\mathrm{CH}^{(2)}_{s}\left(
\mathfrak{Q}^{(2)}
\circ^{0\mathbf{Pth}_{\boldsymbol{\mathcal{A}}^{(2)}}}_{s}
\mathfrak{P}^{(2)}
\right)
\right)$
\allowdisplaybreaks
\begin{align*}
\qquad&=
\mathrm{pr}^{\Theta^{[2]}}_{s}\left(
\mathrm{CH}^{(2)}_{s}\left(
\mathfrak{Q}^{(2)}
\right)
\circ^{0\mathbf{T}_{\Sigma^{\boldsymbol{\mathcal{A}}^{(2)}}}}_{s}
\mathrm{CH}^{(2)}_{s}\left(
\mathfrak{P}^{(2)}
\right)
\right)
\tag{1}
\\&=
\mathrm{pr}^{\Theta^{[2]}}_{s}\left(
\mathrm{CH}^{(2)}_{s}\left(
\mathfrak{Q}^{(2)}
\right)\right)
\circ^{0\mathbf{T}_{\Sigma^{\boldsymbol{\mathcal{A}}^{(2)}}}/{\Theta^{[2]}}}_{s}
\mathrm{pr}^{\Theta^{[2]}}_{s}\left(
\mathrm{CH}^{(2)}_{s}\left(
\mathfrak{P}^{(2)}
\right)\right)
.
\tag{2}
\end{align*}
\end{flushleft}

In the just stated chain of equalities, the first equality follows from the fact that, by Definition~\ref{DDTheta}, the pair
\[
\left(\mathrm{CH}^{(2)}_{s}\left(
\mathfrak{Q}^{(2)}
\circ^{0\mathbf{Pth}_{\boldsymbol{\mathcal{A}}^{(2)}}}_{s}
\mathfrak{P}^{(2)}
\right),
\mathrm{CH}^{(2)}_{s}\left(
\mathfrak{Q}^{(2)}
\right)
\circ^{0\mathbf{T}_{\Sigma^{\boldsymbol{\mathcal{A}}^{(2)}}}(X)}_{s}
\mathrm{CH}^{(2)}_{s}\left(
\mathfrak{P}^{(2)}
\right)
\right)
\]
is in $\Theta^{(2)}_{s}$; finally, the second equality follows from Proposition~\ref{PDThetaCongDCatAlg}.

Hence, $\mathrm{pr}^{\Theta^{[2]}}\circ \mathrm{CH}^{(2)}$ is compatible with the $0$-composition operation.

\textsf{The mapping $\mathrm{pr}^{\Theta^{[2]}}\circ \mathrm{CH}^{(2)}$ is compatible with the second-order rewrite rules.}

Let $s$ be a sort in $S$ and let $\mathfrak{p}^{(2)}$ be a second-order rewrite rule in $\mathcal{A}^{(2)}_{s}$.
The following chain of equalities holds.
\allowdisplaybreaks
\begin{align*}
\mathrm{pr}^{\Theta^{[2]}}_{s}\left(
\mathrm{CH}^{(2)}_{s}\left(
\mathfrak{p}^{(2)\mathbf{Pth}_{\boldsymbol{\mathcal{A}}^{(2)}}}
\right)
\right)&=
\mathrm{pr}^{\Theta^{[2]}}_{s}\left(
\mathfrak{p}^{(2)\mathbf{T}_{\Sigma^{\boldsymbol{\mathcal{A}}^{(2)}}}(X)}
\right)
\tag{1}
\\&=
\mathfrak{p}^{(2)\mathbf{T}_{\Sigma^{\boldsymbol{\mathcal{A}}^{(2)}}}(X)/{\Theta^{[2]}}}.
\tag{2}
\end{align*}

In the just stated chain of equalities, the first equality follows from Proposition~\ref{PDCHDA} and the second equality follows from Proposition~\ref{PDThetaCongDCatAlg}.

Hence, $\mathrm{pr}^{\Theta^{[2]}}\circ \mathrm{CH}^{(2)}$ is compatible with the second-order rewrite rules.

\textsf{The mapping $\mathrm{pr}^{\Theta^{[2]}}\circ \mathrm{CH}^{(2)}$ is compatible with the $1$-source.}

Let $s$ be a sort in $S$ and let us consider the $1$-source operation symbol $\mathrm{sc}^{1}_{s}$ in $\Sigma^{\boldsymbol{\mathcal{A}}^{(2)}}_{s,s}$. Let $\mathfrak{P}^{(2)}$ be a second-order path in $\mathrm{Pth}_{\boldsymbol{\mathcal{A}}^{(2)},s}$. 

The following chain of equalities holds
\allowdisplaybreaks
\begin{align*}
\mathrm{pr}^{\Theta^{[2]}}_{s}\left(
\mathrm{CH}^{(2)}_{s}\left(
\mathrm{sc}^{1\mathbf{Pth}_{\boldsymbol{\mathcal{A}}^{(2)}}}_{s}
\left(
\mathfrak{P}^{(2)}
\right)\right)
\right)&=
\mathrm{pr}^{\Theta^{[2]}}_{s}\left(
\mathrm{sc}^{1\mathbf{T}_{\Sigma^{\boldsymbol{\mathcal{A}}^{(2)}}}}_{s}\left(
\mathrm{CH}^{(2)}_{s}\left(
\mathfrak{P}^{(2)}
\right)\right)
\right)
\tag{1}
\\&=
\mathrm{sc}^{1\mathbf{T}_{\Sigma^{\boldsymbol{\mathcal{A}}^{(2)}}}/{\Theta^{[2]}}}_{s}
\left(
\mathrm{pr}^{\Theta^{[2]}}_{s}\left(
\mathrm{CH}^{(2)}_{s}\left(
\mathfrak{P}^{(2)}
\right)\right)
\right).
\tag{2}
\end{align*}

In the just stated chain of equalities, the first equality follows from the fact that, by Definition~\ref{DDTheta}, the pair
\[
\left(\mathrm{CH}^{(2)}_{s}\left(
\mathrm{sc}^{1\mathbf{Pth}_{\boldsymbol{\mathcal{A}}^{(2)}}}_{s}\left(
\mathfrak{P}^{(2)}
\right)\right),
\mathrm{sc}^{1\mathbf{T}_{\Sigma^{\boldsymbol{\mathcal{A}}^{(2)}}}(X)}_{s}\left(
\mathrm{CH}^{(2)}_{s}\left(
\mathfrak{P}^{(2)}
\right)\right)
\right)
\]
is in $\Theta^{(2)}_{s}$; finally, the second equality follows from Proposition~\ref{PDThetaCongDCatAlg}.

Hence, $\mathrm{pr}^{\Theta^{[2]}}\circ \mathrm{CH}^{(2)}$ is compatible with the $1$-source operation.

\textsf{The mapping $\mathrm{pr}^{\Theta^{[2]}}\circ \mathrm{CH}^{(2)}$ is compatible with the $1$-target.}

Let $s$ be a sort in $S$ and let us consider the $1$-target operation symbol $\mathrm{tg}^{1}_{s}$ in $\Sigma^{\boldsymbol{\mathcal{A}}^{(2)}}_{s,s}$. Let $\mathfrak{P}^{(2)}$ be a second-order path in $\mathrm{Pth}_{\boldsymbol{\mathcal{A}}^{(2)},s}$, then the following equality holds
\[
\mathrm{pr}^{\Theta^{[2]}}_{s}\left(
\mathrm{CH}^{(2)}_{s}\left(
\mathrm{tg}^{1\mathbf{Pth}_{\boldsymbol{\mathcal{A}}^{(2)}}}_{s}\left(
\mathfrak{P}^{(2)}
\right)\right)\right)
=
\mathrm{tg}^{1\mathbf{T}_{\Sigma^{\boldsymbol{\mathcal{A}}^{(2)}}}(X)/{\Theta^{[2]}}}_{s}\left(
\mathrm{pr}^{\Theta^{[2]}}_{s}\left(
\mathrm{CH}^{(2)}_{s}\left(
\mathfrak{P}^{(2)}
\right)\right)\right).
\]

The proof of this case is identical to that of the $1$-source.

Hence, $\mathrm{pr}^{\Theta^{[2]}}\circ \mathrm{CH}^{(2)}$ is compatible with the $1$-target operation.

\textsf{The mapping $\mathrm{pr}^{\Theta^{[2]}}\circ \mathrm{CH}^{(2)}$ is compatible with the $1$-composition.}

Let $s$ be a sort in $S$ and let us consider the $1$-composition operation symbol $\circ^{1}_{s}$ in $\Sigma^{\boldsymbol{\mathcal{A}}^{(2)}}_{ss,s}$. Let $\mathfrak{P}^{(2)}$ and $\mathfrak{Q}^{(2)}$ be two second-order paths in $\mathrm{Pth}_{\boldsymbol{\mathcal{A}}^{(2)},s}$ satisfying that 
\[
\mathrm{sc}^{([1],2)}_{s}\left(\mathfrak{Q}^{(2)}\right)
=
\mathrm{tg}^{([1],2)}_{s}\left(\mathfrak{P}^{(2)}\right).
\]

The following chain of equalities holds
\begin{flushleft}
$\mathrm{pr}^{\Theta^{[2]}}_{s}\left(
\mathrm{CH}^{(2)}_{s}\left(
\mathfrak{Q}^{(2)}
\circ^{1\mathbf{Pth}_{\boldsymbol{\mathcal{A}}^{(2)}}}_{s}
\mathfrak{P}^{(2)}
\right)
\right)$
\allowdisplaybreaks
\begin{align*}
\qquad&=
\mathrm{pr}^{\Theta^{[2]}}_{s}\left(
\mathrm{CH}^{(2)}_{s}\left(
\mathfrak{Q}^{(2)}
\right)
\circ^{1\mathbf{T}_{\Sigma^{\boldsymbol{\mathcal{A}}^{(2)}}}}_{s}
\mathrm{CH}^{(2)}_{s}\left(
\mathfrak{P}^{(2)}
\right)
\right)
\tag{1}
\\&=
\mathrm{pr}^{\Theta^{[2]}}_{s}\left(
\mathrm{CH}^{(2)}_{s}\left(
\mathfrak{Q}^{(2)}
\right)\right)
\circ^{1\mathbf{T}_{\Sigma^{\boldsymbol{\mathcal{A}}^{(2)}}}/{\Theta^{[2]}}}_{s}
\mathrm{pr}^{\Theta^{[2]}}_{s}\left(
\mathrm{CH}^{(2)}_{s}\left(
\mathfrak{P}^{(2)}
\right)\right)
.
\tag{2}
\end{align*}
\end{flushleft}

In the just stated chain of equalities, the first equality follows from the fact that, by Definition~\ref{DDTheta}, the pair
\[
\left(\mathrm{CH}^{(2)}_{s}\left(
\mathfrak{Q}^{(2)}
\circ^{1\mathbf{Pth}_{\boldsymbol{\mathcal{A}}^{(2)}}}_{s}
\mathfrak{P}^{(2)}
\right),
\mathrm{CH}^{(2)}_{s}\left(
\mathfrak{Q}^{(2)}
\right)
\circ^{1\mathbf{T}_{\Sigma^{\boldsymbol{\mathcal{A}}^{(2)}}}(X)}_{s}
\mathrm{CH}^{(2)}_{s}\left(
\mathfrak{P}^{(2)}
\right)
\right)
\]
is in $\Theta^{(2)}_{s}$; finally, the second equality follows from Proposition~\ref{PDThetaCongDCatAlg}.

Hence, $\mathrm{pr}^{\Theta^{[2]}}\circ \mathrm{CH}^{(2)}$ is compatible with the $1$-composition operation.

All in all, we conclude that the composition
$\mathrm{pr}^{\Theta^{[2]}}\circ \mathrm{CH}^{(2)}$
is a $\Sigma^{\boldsymbol{\mathcal{A}}^{(2)}}$-homomorphism from $\mathbf{Pth}_{\boldsymbol{\mathcal{A}}^{(2)}}$ to $\mathbf{T}_{\Sigma^{\boldsymbol{\mathcal{A}}^{(2)}}}(X)/{\Theta^{[2]}}$.

This completes the proof.
\end{proof}

The following lemma states that if a term of sort $s$ is such that its image under $\mathrm{ip}^{(2,X)@}_{s}$ is a second-order path, then it is related, with respect to the $\Sigma^{\boldsymbol{\mathcal{A}}^{(2)}}$-congruence  $\Theta^{[2]}$ on $\mathbf{T}_{\Sigma^{\boldsymbol{\mathcal{A}}^{(2)}}}(X)$, with a term of the same sort in $\mathrm{CH}^{(2)}[\mathrm{Pth}_{\boldsymbol{\mathcal{A}}^{(2)}}]$. Actually, we prove that such a term is related with its image under the action of $\mathrm{CH}^{(2)}_{s}\circ\mathrm{ip}^{(2,X)@}_{s}$.

\begin{restatable}{lemma}{LDWCong}
\label{LDWCong} Let $s$ be a sort in $S$ and $P$ a term in $\mathrm{T}_{\Sigma^{\boldsymbol{\mathcal{A}}^{(2)}}}(X)_{s}$. If $\mathrm{ip}^{(2,X)@}_{s}(P)$ is a second-order path in $\mathrm{Pth}_{\boldsymbol{\mathcal{A}}^{(2)},s}$, then $(P,\mathrm{CH}^{(2)}_{s}(\mathrm{ip}^{(2,X)@}_{s}(P)))\in\Theta^{[2]}_{s}$. 
\end{restatable}
\begin{proof}
We prove it by induction on the height of $P$.

\textsf{Base step of the induction.}

If $z\in\mathrm{B}^{0}_{\Sigma^{\boldsymbol{\mathcal{A}}^{(2)}}}(X)_{s}$, then either (1) $z\in\eta^{(2,0)\sharp}[\mathrm{B}^{0}_{\Sigma}(X)]_{s}$, or~(2) $z=\mathfrak{p}^{\mathbf{T}_{\Sigma^{\boldsymbol{\mathcal{A}}^{(2)}}}(X)}$, for some rewrite rule $\mathfrak{p}\in\mathcal{A}_{s}$, or~(3) $z=\mathfrak{p}^{(2)\mathbf{T}_{\Sigma^{\boldsymbol{\mathcal{A}}^{(2)}}}(X)}$, for some second-order rewrite rule $\mathfrak{p}^{(2)}\in\mathcal{A}^{(2)}_{s}$. Let us recall that in either case, $\mathrm{ip}^{(2,X)@}_{s}(z)$ is a second-order path in $\mathrm{Pth}_{\boldsymbol{\mathcal{A}}^{(2)},s}$.

If~(1), then $z=\eta^{(2,0)\sharp}_{s}(x)$ for some term $x\in\mathrm{B}^{0}_{\Sigma}(X)_{s}$.

Hence the following chain of equalities holds
\allowdisplaybreaks
\begin{align*}
\mathrm{CH}^{(2)}_{s}\left(
\mathrm{ip}^{(2,X)@}_{s}\left(
z\right)\right)
&=
\mathrm{CH}^{(2)}_{s}\left(\mathrm{ip}^{(2,[1])\sharp}_{s}\left(
[\eta^{(1,0)\sharp}_{s}(x)]_{s}
\right)
\right)
\tag{1}
\\&=
\eta^{(2,1)\sharp}_{s}\left(
\mathrm{CH}^{(1)\mathrm{m}}_{s}\left(
\mathrm{ip}^{([1],X)@}_{s}\left(
[\eta^{(1,0)\sharp}_{s}(x)]_{s}
\right)\right)\right)
\tag{2}
\\&=
\eta^{(2,1)\sharp}_{s}\left(
\mathrm{CH}^{(1)\mathrm{m}}_{s}\left(
\left[\mathrm{ip}^{(1,X)@}_{s}\left(
\eta^{(1,0)\sharp}_{s}(x)
\right)
\right]_{s}
\right)\right)
\tag{3}
\\&=
\eta^{(2,1)\sharp}_{s}\left(
\mathrm{CH}^{(1)\mathrm{m}}_{s}\left(
\left[\mathrm{ip}^{(1,0)\sharp}_{s}(x)
\right]_{s}
\right)\right)
\tag{4}
\\&=
\eta^{(2,1)\sharp}_{s}\left(
\mathrm{CH}^{(1)}_{s}\left(
\mathrm{ip}^{(1,0)\sharp}_{s}(x)
\right)\right)
\tag{5}
\\&=
\eta^{(2,1)\sharp}_{s}\left(
\eta^{(1,0)\sharp}_{s}\left(
x
\right)\right)
\tag{6}
\\&=
\eta^{(2,0)\sharp}_{s}(x)
\tag{7}
\\&=
z.
\tag{8}
\end{align*}

The first equality follows from the fact that the image of $\eta^{(2,0)\sharp}_{s}(x)$ under the  $\mathrm{ip}^{(2,X)@}$ mapping is, by Definition~\ref{DDIp} and Definition~\ref{DDPth}, the $(2,[1])$-identity second-order path on $[\eta^{(1,0)\sharp}_{s}(x)]_{s}$; the second equality follows from Proposition~\ref{PDCHDUId}; the third equality applies the mapping $\mathrm{ip}^{([1],X)@}$ introduced in Definition~\ref{DPTQIp}; the fourth equality follows from Proposition~\ref{PIpUZ}; the fifth equality applies the monomorphic Curry-Howard mapping at a path class; the sixth equality follows from Proposition~\ref{PCHId}; the seventh equality follows from Proposition~\ref{PDEmb}; finally, the last equality recovers the definition of $z$.

By reflexivity, we conclude that $(z,\mathrm{CH}^{(2)}_{s}(\mathrm{ip}^{(2,X)@}_{s}(z)))\in\Theta^{[2]}_{s}$.

\label{********ICI}
If~(2), then  $z=\mathfrak{p}^{\mathbf{T}_{\Sigma^{\boldsymbol{\mathcal{A}}^{(2)}}}(X)}$, for some rewrite rule $\mathfrak{p}\in\mathcal{A}_{s}$. 

Hence the following chain of equalities holds
\allowdisplaybreaks
\begin{align*}
\mathrm{CH}^{(2)}_{s}\left(
\mathrm{ip}^{(2,X)@}_{s}\left(
\mathfrak{p}^{\mathbf{T}_{\Sigma^{\boldsymbol{\mathcal{A}}^{(2)}}}(X)}
\right)\right)
&=
\mathrm{CH}^{(2)}_{s}\left(
\mathfrak{p}^{\mathbf{F}_{\Sigma^{\boldsymbol{\mathcal{A}}^{(2)}}}\left(
\mathbf{Pth}_{\boldsymbol{\mathcal{A}}^{(2)}}
\right)}\right)
\tag{1}
\\&=
\mathrm{CH}^{(2)}_{s}\left(
\mathfrak{p}^{\mathbf{Pth}_{\boldsymbol{\mathcal{A}}^{(2)}}}
\right)
\tag{2}
\\&=
\mathrm{CH}^{(2)}_{s}\left(
\mathrm{ip}^{(2,[1])\sharp}_{s}\left(
\left[\mathfrak{p}^{\mathbf{PT}_{\boldsymbol{\mathcal{A}}}}
\right]
_{s}
\right)\right)
\tag{3}
\\&=
\eta^{(2,1)\sharp}_{s}\left(
\mathrm{CH}^{(1)\mathrm{m}}_{s}\left(
\mathrm{ip}^{([1],X)@}_{s}\left(
\left[
\mathfrak{p}^{\mathbf{PT}_{\boldsymbol{\mathcal{A}}}}
\right]_{s}
\right)\right)\right)
\tag{4}
\\&=
\eta^{(2,1)\sharp}_{s}\left(
\mathrm{CH}^{(1)\mathrm{m}}_{s}\left(
\left[
\mathrm{ip}^{(1,X)@}_{s}\left(
\mathfrak{p}^{\mathbf{PT}_{\boldsymbol{\mathcal{A}}}}
\right)\right]_{s}
\right)\right)
\tag{5}
\\&=
\eta^{(2,1)\sharp}_{s}\left(
\mathrm{CH}^{(1)\mathrm{m}}_{s}\left(
\left[\mathfrak{p}^{
\mathbf{Pth}_{\boldsymbol{\mathcal{A}}}
}
\right]_{s}
\right)\right)
\tag{6}
\\&=
\eta^{(2,1)\sharp}_{s}\left(
\mathfrak{p}^{
\mathbf{T}_{\Sigma^{\boldsymbol{\mathcal{A}}}}(X)
}
\right)
\tag{7}
\\&=
\mathfrak{p}^{\mathbf{T}_{\Sigma^{\boldsymbol{\mathcal{A}}^{(2)}}}(X)}.
\tag{8}
\end{align*}

The first equality follows from the fact that, by Definition~\ref{DDIp}, $\mathrm{ip}^{(2,X)@}$ is a $\Sigma^{\boldsymbol{\mathcal{A}}^{(2)}}$-homomorphism; the second equality follows from the fact that the operation symbol $\mathfrak{p}$ in $\Sigma^{\boldsymbol{\mathcal{A}}^{(2)}}_{\lambda,s}$ is defined in the many-sorted partial $\Sigma^{\boldsymbol{\mathcal{A}}^{(2)}}$-algebra $\mathbf{Pth}_{\boldsymbol{\mathcal{A}}^{(2)}}$ according to Proposition~\ref{PDPthDCatAlg}, thus the interpretation of the operation symbol $\mathfrak{p}$ in the free completion $\mathbf{F}_{\Sigma^{\boldsymbol{\mathcal{A}}^{(2)}}}(\mathbf{Pth}_{\boldsymbol{\mathcal{A}}^{(2)}})$ becomes that interpretation occurring in $\mathbf{Pth}_{\boldsymbol{\mathcal{A}}^{(2)}}$; the third equality unravels the interpretation of the operation symbol $\mathfrak{p}$ in $\mathbf{Pth}_{\boldsymbol{\mathcal{A}}^{(2)}}$, according to Proposition~\ref{PDPthDCatAlg}. Let us recall that $\mathfrak{p}^{\mathbf{Pth}_{\boldsymbol{\mathcal{A}}^{(2)}}}$ is given by the $(2,[1])$-identity second-order path on the path term class $[\mathfrak{p}^{\mathbf{PT}_{\boldsymbol{\mathcal{A}}}}]_{s}$; the fourth equality follows from Proposition~\ref{PDCHDUId}; the fifth equality applies the mapping $\mathrm{ip}^{([1],X)@}$ introduced in Definition~\ref{DPTQIp}. Note that, by Definition~\ref{DIp}, $\mathrm{ip}^{(1,X)@}$ is a $\Sigma^{\boldsymbol{\mathcal{A}}}$-homomorphism. Let us recall that the operation symbol $\mathfrak{p}$ in $\Sigma^{\boldsymbol{\mathcal{A}}}_{\lambda,s}$ is defined in the many-sorted partial $\Sigma^{\boldsymbol{\mathcal{A}}}$-algebra $\mathbf{Pth}_{\boldsymbol{\mathcal{A}}}$ according to Proposition~\ref{PPthCatAlg}; the seventh equality applies the monomorphic Curry-Howard mapping at a path class. Let us recall that the resulting term is a consequence of Proposition~\ref{PCHA}; finally, the last equality follows from the fact that $\eta^{(2,1)\sharp}$ is a $\Sigma^{\boldsymbol{\mathcal{A}}}$-homomorphism according to Proposition~\ref{PDUEmb}.

By reflexivity, $(\mathfrak{p}^{\mathbf{T}_{\Sigma^{\boldsymbol{\mathcal{A}}^{(2)}}}(X)},\mathrm{CH}^{(2)}_{s}(\mathrm{ip}^{(2,X)@}_{s}(\mathfrak{p}^{\mathbf{T}_{\Sigma^{\boldsymbol{\mathcal{A}}^{(2)}}}(X)})))\in\Theta^{[2]}_{s}$.

If~(3), then  $z=\mathfrak{p}^{(2)\mathbf{T}_{\Sigma^{\boldsymbol{\mathcal{A}}^{(2)}}}(X)}$, for some second-order rewrite rule $\mathfrak{p}^{(2)}\in\mathcal{A}^{(2)}_{s}$. 

Hence the following chain of equalities holds
\allowdisplaybreaks
\begin{align*}
\mathrm{CH}^{(2)}_{s}\left(
\mathrm{ip}^{(2,X)@}_{s}\left(
\mathfrak{p}^{(2)\mathbf{T}_{\Sigma^{\boldsymbol{\mathcal{A}}^{(2)}}}(X)}
\right)\right)
&=
\mathrm{CH}^{(2)}_{s}\left(
\mathfrak{p}^{(2)\mathbf{F}_{\Sigma^{\boldsymbol{\mathcal{A}}^{(2)}}}
\left(\mathbf{Pth}_{\boldsymbol{\mathcal{A}}^{(2)}}
\right)
}\right)
\tag{1}
\\&=
\mathrm{CH}^{(2)}_{s}\left(
\mathfrak{p}^{(2)\mathbf{Pth}_{\boldsymbol{\mathcal{A}}^{(2)}}}
\right)
\tag{2}
\\&=
\mathrm{CH}^{(2)}_{s}\left(
\mathrm{ech}^{(2,\mathcal{A}^{(2)})}_{s}\left(\mathfrak{p}^{(2)}
\right)\right)
\tag{3}
\\&=
\mathfrak{p}^{(2)\mathbf{T}_{\Sigma^{\boldsymbol{\mathcal{A}}^{(2)}}}(X)}.
\tag{4}
\end{align*}

The first equality follows from the fact that $\mathrm{ip}^{(2,X)@}$ is a $\Sigma^{\boldsymbol{\mathcal{A}}^{(2)}}$-homomorphism in virtue of Definition~\ref{DDIp}; the second equality follows from the fact that the operation symbol $\mathfrak{p}^{(2)}$ in $\Sigma^{\boldsymbol{\mathcal{A}}^{(2)}}_{\lambda,s}$ is defined in the many-sorted partial $\Sigma^{\boldsymbol{\mathcal{A}}^{(2)}}$-algebra $\mathbf{Pth}_{\boldsymbol{\mathcal{A}}^{(2)}}$ according to Proposition~\ref{PDPthDCatAlg}, thus the interpretation of the operation symbol $\mathfrak{p}^{(2)}$ in the free completion $\mathbf{F}_{\Sigma^{\boldsymbol{\mathcal{A}}^{(2)}}}(\mathbf{Pth}_{\boldsymbol{\mathcal{A}}^{(2)}})$ becomes that interpretation occurring in $\mathbf{Pth}_{\boldsymbol{\mathcal{A}}^{(2)}}$; the third equality recovers the interpretation of the operation symbol $\mathfrak{p}^{(2)}$ in $\mathbf{Pth}_{\boldsymbol{\mathcal{A}}^{(2)}}$ according to Proposition~\ref{PDPthDCatAlg}; finally, the last equality follows from Proposition~\ref{PDCHA}.

By reflexivity, $(\mathfrak{p}^{(2)\mathbf{T}_{\Sigma^{\boldsymbol{\mathcal{A}}^{(2)}}}(X)},\mathrm{CH}^{(2)}_{s}(\mathrm{ip}^{(2,X)@}_{s}(\mathfrak{p}^{(2)\mathbf{T}_{\Sigma^{\boldsymbol{\mathcal{A}}^{(2)}}}(X)})))\in\Theta^{[2]}_{s}$.

\textsf{Inductive step of the induction.}

Assume that the statement holds for terms up to height $n\in\mathbb{N}$, i.e., for every term $P\in\mathrm{T}_{\Sigma^{\boldsymbol{\mathcal{A}}^{(2)}}}(X)_{s}$ with $P\in\mathrm{B}^{n}_{\Sigma^{\boldsymbol{\mathcal{A}}^{(2)}}}(X)_{s}$, if $\mathrm{ip}^{(2,X)@}_{s}(P)$ is a second-order path in $\mathrm{Pth}_{\boldsymbol{\mathcal{A}}^{(2)},s}$, then $(P,\mathrm{CH}^{(2)}_{s}(\mathrm{ip}^{(2,X)@}_{s}(P)))\in\Theta^{[2]}_{s}$.

Now, we prove it for the case in which $P$ is a term in $\mathrm{T}_{\Sigma^{\boldsymbol{\mathcal{A}}^{(2)}}}(X)_{s}$ with $P\in\mathrm{B}^{n+1}_{\Sigma^{\boldsymbol{\mathcal{A}}^{(2)}}}(X)_{s}$ satisfying that $\mathrm{ip}^{(2,X)@}_{s}(P)$ is a second-order path in $\mathrm{Pth}_{\boldsymbol{\mathcal{A}}^{(2)},s}$. Since $P\in\mathrm{B}^{n+1}_{\Sigma^{\boldsymbol{\mathcal{A}}^{(2)}}}(X)_{s}$, then there exists a unique $\mathbf{s}\in S^{\star}-\{\lambda\}$, a unique operation $\gamma\in\Sigma^{\boldsymbol{\mathcal{A}}^{(2)}}_{\mathbf{s},s}$ and a unique family of terms $(P_{j})_{j\in\bb{\mathbf{s}}}\in\mathrm{B}^{n}_{\Sigma^{\boldsymbol{\mathcal{A}}^{(2)}}}(X)_{\mathbf{s}}$ satisfying that 
$$
P=\gamma^{\mathbf{T}_{\Sigma^{\boldsymbol{\mathcal{A}}^{(2)}}}(X)}
\left(\left(
P_{j}
\right)_{j\in\bb{\mathbf{s}}}\right).
$$

Since we are assuming that $\mathrm{ip}^{(2,X)@}_{s}(P)$ is a second-order path in $\mathrm{Pth}_{\boldsymbol{\mathcal{A}}^{(2)},s}$ and $P$ is equal to $\gamma^{\mathbf{T}_{\Sigma^{\boldsymbol{\mathcal{A}}^{(2)}}}(X)}
((
P_{j}
)_{j\in\bb{\mathbf{s}}})$, we conclude, by Lemma~\ref{LDThetaCongSub}, that, for every $j\in\bb{\mathbf{s}}$, $\mathrm{ip}^{(2,X)@}_{s_{j}}(P_{j})$ is a second-order path in $\mathrm{Pth}_{\boldsymbol{\mathcal{A}}^{(2)},s_{j}}$. Moreover, by induction, we have that
$$
\left(P_{j},
\mathrm{CH}^{(2)}_{s_{j}}\left(
\mathrm{ip}^{(2,X)@}_{s_{j}}\left(
P_{j}
\right)\right)
\right)\in\Theta^{[2]}_{s_{j}}.
$$

Taking into account that $\Theta^{[2]}$ is a $\Sigma^{\boldsymbol{\mathcal{A}}^{(2)}}$-congruence on $\mathbf{T}_{\Sigma^{\boldsymbol{\mathcal{A}}^{(2)}}}(X)$, then
$$
\left(P,
\gamma^{\mathbf{T}_{\Sigma^{\boldsymbol{\mathcal{A}}^{(2)}}}(X)}
\left(\left(
\mathrm{CH}^{(2)}_{s_{j}}\left(
\mathrm{ip}^{(2,X)@}_{s_{j}}\left(
P_{j}
\right)\right)\right)_{j\in\bb{\mathbf{s}}}\right)
\right)\in\Theta^{[2]}_{s}.
$$

On the other hand, note that the following chain of equalities hold
\begin{flushleft}
$\mathrm{pr}^{\Theta^{[2]}}_{s}\left(
\mathrm{CH}^{(2)}_{s}\left(
\mathrm{ip}^{(2,X)@}_{s}\left(
P
\right)
\right)
\right)
$
\allowdisplaybreaks
\begin{align*}
\qquad&=
\mathrm{pr}^{\Theta^{[2]}}_{s}\left(
\mathrm{CH}^{(2)}_{s}\left(
\mathrm{ip}^{(2,X)@}_{s}\left(
\gamma^{\mathbf{T}_{\Sigma^{\boldsymbol{\mathcal{A}}^{(2)}}}(X)}
\left(
\left(
P_{j}
\right)_{j\in\bb{\mathbf{s}}}
\right)
\right)
\right)
\right)
\tag{1}
\\&=
\mathrm{pr}^{\Theta^{[2]}}_{s}\left(
\mathrm{CH}^{(2)}_{s}\left(
\gamma^{\mathbf{Pth}_{\boldsymbol{\mathcal{A}}^{(2)}}}
\left(
\left(
\mathrm{ip}^{(2,X)@}_{s_{j}}\left(
P_{j}
\right)
\right)_{j\in\bb{\mathbf{s}}}
\right)
\right)
\right)
\tag{2}
\\&=
\gamma^{\mathbf{T}_{\Sigma^{\boldsymbol{\mathcal{A}}^{(2)}}}(X)/{\Theta^{[2]}}}
\left(
\left(
\mathrm{pr}^{\Theta^{[2]}}_{s_{j}}\left(
\mathrm{CH}^{(2)}_{s_{j}}\left(
\mathrm{ip}^{(2,X)@}_{s_{j}}\left(
P_{j}
\right)
\right)
\right)
\right)_{j\in\bb{\mathbf{s}}}
\right)
\tag{3}
\\&=
\mathrm{pr}^{\Theta^{[2]}}_{s}\left(
\gamma^{\mathbf{T}_{\Sigma^{\boldsymbol{\mathcal{A}}^{(2)}}}(X)}
\left(
\left(
\mathrm{CH}^{(2)}_{s_{j}}\left(
\mathrm{ip}^{(2,X)@}_{s_{j}}\left(
P_{j}
\right)
\right)
\right)_{j\in\bb{\mathbf{s}}}
\right)
\right).
\tag{4}
\end{align*}
\end{flushleft}

In the just stated chain of equalities, the first equality unravels the description of $P$; the second equality follows from the fact that, by assumption $\mathrm{ip}^{(2,X)@}_{s}(P)$ is a second-order path and, by Definition~\ref{DDIp}, $\mathrm{ip}^{(2,X)@}$ is a $\Sigma^{\boldsymbol{\mathcal{A}}^{(2)}}$-homomorphism; the third equality follows from the fact that, according to Proposition~\ref{PDThetaDCH}, $\mathrm{pr}^{\Theta^{[2]}}\circ \mathrm{CH}^{(2)}$ is a $\Sigma^{\boldsymbol{\mathcal{A}}^{(2)}}$-homomorphism; finally, the last equality follows from the fact that, according to Proposition~\ref{PDThetaCongDCatAlg}, $\mathrm{pr}^{\Theta^{[2]}}$ is a $\Sigma^{\boldsymbol{\mathcal{A}}^{(2)}}$-homomorphism.

Consequently, the following pair is in $\Theta^{[2]}$
\[
\left(
\mathrm{CH}^{(2)}_{s}\left(
\mathrm{ip}^{(2,X)@}_{s}\left(
P
\right)
\right)
,
\gamma^{\mathbf{T}_{\Sigma^{\boldsymbol{\mathcal{A}}^{(2)}}}(X)}
\left(
\left(
\mathrm{CH}^{(2)}_{s_{j}}\left(
\mathrm{ip}^{(2,X)@}_{s_{j}}\left(
P_{j}
\right)
\right)
\right)_{j\in\bb{\mathbf{s}}}
\right)
\right)
\in\Theta^{[2]}_{s}.
\]

Since $\Theta^{[2]}_{s}$ is an equivalence relation, we conclude that 
\[
\left(
P,
\mathrm{CH}^{(2)}_{s}\left(
\mathrm{ip}^{(2,X)@}_{s}\left(
P
\right)
\right)
\right)\in\Theta^{[2]}_{s}.
\]

This completes the proof.
\end{proof}

We continue with the description of the behaviour of the $\Sigma^{\boldsymbol{\mathcal{A}}^{(2)}}$-congruence $\Theta^{[2]}$. In the following lemma we state that, for every sort $s\in S$, if two terms are $\Theta^{[2]}_{s}$-related and one of them, when mapped under $\mathrm{ip}^{(2,X)@}_{s}$, retrieves a second-order path, then the other term has a similar behaviour. Moreover, if the just described situation happens, then these two second-order paths will have the same image under the second-order Curry-Howard mapping.

\begin{restatable}{lemma}{LDThetaCong}
\label{LDThetaCong} Let $s$ be a sort in $S$ and $P,Q$ terms in $\mathrm{T}_{\Sigma^{\boldsymbol{\mathcal{A}}^{(2)}}}(X)_{s}$ such that $(P,Q)\in\Theta^{[2]}_{s}$, then
\begin{itemize}
\item[(i)] $\mathrm{ip}^{(2,X)@}_{s}(P)\in\mathrm{Pth}_{\boldsymbol{\mathcal{A}}^{(2)},s}$ if, and only if, $\mathrm{ip}^{(2,X)@}_{s}(Q)\in\mathrm{Pth}_{\boldsymbol{\mathcal{A}}^{(2)},s}$;
\item[(ii)] If $\mathrm{ip}^{(2,X)@}_{s}(P)\in\mathrm{Pth}_{\boldsymbol{\mathcal{A}}^{(2)},s}$ or  $\mathrm{ip}^{(2,X)@}_{s}(Q)\in\mathrm{Pth}_{\boldsymbol{\mathcal{A}}^{(2)},s}$ is a second-order path in $\mathrm{Pth}_{\boldsymbol{\mathcal{A}}^{(2)},s}$ then
$
\mathrm{CH}^{(2)}_{s}(
\mathrm{ip}^{(2,X)@}_{s}(
P))=\mathrm{CH}^{(2)}_{s}(
\mathrm{ip}^{(2,X)@}_{s}(
Q)).
$
\end{itemize}
\end{restatable}
\begin{proof}
We recall from Remark~\ref{RDThetaCong} that 
$$
\Theta^{[2]}=
\bigcup_{n\in\mathbb{N}}\mathrm{C}^{n}_{\Sigma^{\boldsymbol{\mathcal{A}}^{(2)}}}
\left(
\Theta^{(2)}\right).
$$

We prove the statement by induction on $n\in\mathbb{N}$.

\textsf{Base step of the induction.}

Let us recall from Remark~\ref{RDThetaCong} that 
\[
\mathrm{C}^{0}_{\Sigma^{\boldsymbol{\mathcal{A}}^{(2)}}}\left(\Theta^{(2)}\right)=
\left(\Theta^{(2)}\right)
\cup
\left(
\Theta^{(2)}
\right)^{-1}
\cup
\Delta_{\mathrm{T}_{\Sigma^{\boldsymbol{\mathcal{A}}^{(2)}}}(X)}.
\]

The statement trivially holds for a pair $(P,Q)$ in $\Delta_{\mathrm{T}_{\Sigma^{\boldsymbol{\mathcal{A}}^{(2)}}}(X)_{s}}$.

If the pair $(P,Q)$ is in $\Theta^{(2)}_{s}$, following Definition~\ref{DDTheta}, in all cases, the term $P$ belongs to $\mathrm{CH}^{(2)}_{s}[\mathrm{Pth}_{\boldsymbol{\mathcal{A}}^{(2)},s}]$. Following Proposition~\ref{PDIpDCH}, we know that $\mathrm{ip}^{(2,X)@}_{s}(P)$ is a second-order path in $\mathrm{Pth}_{\boldsymbol{\mathcal{A}}^{(2)},s}$. The different cases for $Q$ are handled by Lemmas~\ref{LDThetaScZ},~\ref{LDThetaTgZ},~\ref{LDThetaCompZ},~\ref{LDThetaScU},~\ref{LDThetaTgU} or~\ref{LDThetaCompU}. In any case, it follows that $\mathrm{ip}^{(2,X)@}_{s}(Q)$ is also a second-order path in $\mathrm{Pth}_{\boldsymbol{\mathcal{A}}^{(2)},s}$. Moreover,
$
\mathrm{CH}^{(2)}_{s}(
\mathrm{ip}^{(2,X)@}_{s}(
P))=
\mathrm{CH}^{(2)}_{s}(
\mathrm{ip}^{(2,X)@}_{s}(
Q)).
$ 

In case $(P,Q)$ is a pair in $\left(\Theta^{(2)}\right)^{-1}_{s}$ we reason in a similar way.

This completes the base case.

\textsf{Inductive step of the induction.}

Assume the statement holds for $n\in\mathbb{N}$, i.e., for every sort $s\in S$ and every pair of terms $(P,Q)$ in $\mathrm{T}_{\Sigma^{\boldsymbol{\mathcal{A}}^{(2)}}}(X)_{s}$ such that $(P,Q)$ in $\mathrm{C}^{n}_{\Sigma^{\boldsymbol{\mathcal{A}}^{(2)}}}(\Theta^{(2)})_{s}$ then
\begin{itemize}
\item[(i)] $\mathrm{ip}^{(2,X)@}_{s}(P)\in\mathrm{Pth}_{\boldsymbol{\mathcal{A}}^{(2)},s}$ if, and only if, $\mathrm{ip}^{(2,X)@}_{s}(Q)\in\mathrm{Pth}_{\boldsymbol{\mathcal{A}}^{(2)},s}$;
\item[(ii)] If $\mathrm{ip}^{(2,X)@}_{s}(P)\in\mathrm{Pth}_{\boldsymbol{\mathcal{A}}^{(2)},s}$ or  $\mathrm{ip}^{(2,X)@}_{s}(Q)\in\mathrm{Pth}_{\boldsymbol{\mathcal{A}}^{(2)},s}$ is a second-order path in $\mathrm{Pth}_{\boldsymbol{\mathcal{A}}^{(2)},s}$ then
$
\mathrm{CH}^{(2)}_{s}(
\mathrm{ip}^{(2,X)@}_{s}(
P))=\mathrm{CH}^{(2)}_{s}(
\mathrm{ip}^{(2,X)@}_{s}(
Q)).
$
\end{itemize}

We prove the statement for $n+1$. Let $s$ be a sort in $S$ and let $(P,Q)$ be a pair of terms in $\mathrm{T}_{\Sigma^{\boldsymbol{\mathcal{A}}^{(2)}}}(X)_{s}$ such that $(P,Q)$ is a pair in $\mathrm{C}^{n+1}_{\Sigma^{\boldsymbol{\mathcal{A}}^{(2)}}}(\Theta^{(2)})_{s}$. Let us recall from Remark~\ref{RDThetaCong} that 
\begin{multline*}
\mathrm{C}^{n+1}_{\Sigma^{\boldsymbol{\mathcal{A}}^{(2)}}}\left(
\Theta^{(2)}\right)_{s}=
\left(\mathrm{C}^{n}_{\Sigma^{\boldsymbol{\mathcal{A}}^{(2)}}}\left(
\Theta^{(2)}
\right)_{s}
\circ
\mathrm{C}^{n}_{\Sigma^{\boldsymbol{\mathcal{A}}^{(2)}}}\left(
\Theta^{(2)}
\right)_{s}
\right)\cup
\\
\left(
\bigcup_{\gamma\in\Sigma^{\boldsymbol{\mathcal{A}}^{(2)}}_{\neq\lambda,s}}
\gamma^{\mathbf{T}_{\Sigma^{\boldsymbol{\mathcal{A}}^{(2)}}}(X)}
\times
\gamma^{\mathbf{T}_{\Sigma^{\boldsymbol{\mathcal{A}}^{(2)}}}(X)}
\left[
\mathrm{C}^{n}_{\Sigma^{\boldsymbol{\mathcal{A}}^{(2)}}}\left(\Theta^{(2)}\right)_{
\mathrm{ar}(\gamma)
}
\right]
\right)_{s\in S}.
\end{multline*}

Then either~(1), $(P,Q)$ is in $\mathrm{C}^{n}_{\Sigma^{\boldsymbol{\mathcal{A}}^{(2)}}}(\Theta^{(2)})_{s}
\circ
\mathrm{C}^{n}_{\Sigma^{\boldsymbol{\mathcal{A}}^{(2)}}}(\Theta^{(2)})_{s}
$ or~(2), $(P,Q)$ is in $\gamma^{\mathbf{T}_{\Sigma^{\boldsymbol{\mathcal{A}}^{(2)}}}(X)}
\times
\gamma^{\mathbf{T}_{\Sigma^{\boldsymbol{\mathcal{A}}^{(2)}}}(X)}
\left[
\mathrm{C}^{n}_{\Sigma^{\boldsymbol{\mathcal{A}}^{(2)}}}\left(\Theta^{(2)}\right)_{
\mathrm{ar}(\gamma)
}
\right]$ for some operation symbol $\gamma\in\Sigma^{\boldsymbol{\mathcal{A}}^{(2)}}_{\neq\lambda,s}$.

If~(1), then there exists a term $R\in\mathrm{T}_{\Sigma^{\boldsymbol{\mathcal{A}}^{(2)}}}(X)_{s}$ for which $(P,R)$ and $(R,Q)$ belong to $\mathrm{C}^{n}_{\Sigma^{\boldsymbol{\mathcal{A}}^{(2)}}}(\Theta^{(2)})_{s}$. Then by induction we have that
$$
\mathrm{ip}^{(2,X)@}_{s}(P)
\in\mathrm{Pth}_{\boldsymbol{\mathcal{A}}^{(2)},s}
\Leftrightarrow
\mathrm{ip}^{(2,X)@}_{s}(R)
\in\mathrm{Pth}_{\boldsymbol{\mathcal{A}}^{(2)},s}
\Leftrightarrow
\mathrm{ip}^{(2,X)@}_{s}(Q)
\in\mathrm{Pth}_{\boldsymbol{\mathcal{A}}^{(2)},s}.
$$

Moreover, in case one of the previous elements is a second-order path in $\mathrm{Pth}_{\boldsymbol{\mathcal{A}}^{(2)},s}$, we have, by induction, that
$$
\mathrm{CH}^{(2)}_{s}\left(
\mathrm{ip}^{(2,X)@}_{s}\left(
P
\right)\right)
=
\mathrm{CH}^{(2)}_{s}\left(
\mathrm{ip}^{(2,X)@}_{s}\left(
R
\right)\right)
=
\mathrm{CH}^{(2)}_{s}\left(
\mathrm{ip}^{(2,X)@}_{s}\left(
Q
\right)\right).
$$

If~(2), then there exists a unique word $\mathbf{s}\in S^{\star}-\{\lambda\}$, a unique operation symbol $\gamma\in \Sigma^{\boldsymbol{\mathcal{A}}^{(2)}}_{\mathbf{s},s}$ and a unique family of pairs $((P_{j},Q_{j}))_{j\in\bb{\mathbf{s}}}$ in $\mathrm{C}^{n}_{\Sigma^{\boldsymbol{\mathcal{A}}^{(2)}}}(\Theta^{(2)})_{\mathbf{s}}$ for which
$$
(P,Q)=
\left(\gamma^{\mathbf{T}_{\Sigma^{\boldsymbol{\mathcal{A}}^{(2)}}}(X)}
\left(\left(
P_{j}
\right)_{j\in\bb{\mathbf{s}}}\right),
\gamma^{\mathbf{T}_{\Sigma^{\boldsymbol{\mathcal{A}}^{(2)}}}(X)}
\left(\left(
Q_{j}
\right)_{j\in\bb{\mathbf{s}}}\right)
\right).
$$

We will distinguish the following cases according to the different possibilities for the operation symbol $\gamma\in\Sigma^{\boldsymbol{\mathcal{A}}^{(2)}}_{\mathbf{s},s}$. Note that either $\gamma$ is an operation symbol $\sigma\in\Sigma_{\mathbf{s},s}$ or the operation symbol of $0$-source, $\mathrm{sc}^{0}_{s}$, or the operation symbol of $0$-target, $\mathrm{tg}^{0}_{s}$, or the operation symbol of $0$-composition, $\circ^{0}_{s}$, or the operation symbol of $1$-source, $\mathrm{sc}^{1}_{s}$, or the operation symbol of $1$-target, $\mathrm{tg}^{1}_{s}$, or  the operation symbol of $1$-composition $\circ^{1}_{s}$.

\textsf{$\gamma$ is an operation symbol $\sigma\in\Sigma_{\mathbf{s},s}$.}

Let $\mathbf{s}$ be a word in $S^{\star}-\{\lambda\}$, let $\sigma$ be an operation symbol in $\Sigma_{\mathbf{s},s}$ and let $((P_{j},Q_{j}))_{j\in\bb{\mathbf{s}}}$ be a family of pairs in $\mathrm{C}^{n}_{\Sigma^{\boldsymbol{\mathcal{A}}^{(2)}}}(\Theta^{(2)})_{\mathbf{s}}$ for which 
$$
(P,Q)=
\left(\sigma^{\mathbf{T}_{\Sigma^{\boldsymbol{\mathcal{A}}^{(2)}}}(X)}
\left(\left(
P_{j}
\right)_{j\in\bb{\mathbf{s}}}\right),
\sigma^{\mathbf{T}_{\Sigma^{\boldsymbol{\mathcal{A}}^{(2)}}}(X)}
\left(\left(
Q_{j}
\right)_{j\in\bb{\mathbf{s}}}\right)
\right).
$$

Note that the following chain of equivalences holds
\allowdisplaybreaks
\begin{align*}
\mathrm{ip}^{(2,X)@}_{s}(P)\in\mathrm{Pth}_{\boldsymbol{\mathcal{A}}^{(2)},s}
&\Leftrightarrow
\mbox{for every }j\in\bb{\mathbf{s}},\,\mathrm{ip}^{(2,X)@}_{s_{j}}(P_{j})\in\mathrm{Pth}_{\boldsymbol{\mathcal{A}}^{(2)},s_{j}}
\tag{1}
\\&\Leftrightarrow
\mbox{for every }j\in\bb{\mathbf{s}},\,\mathrm{ip}^{(2,X)@}_{s_{j}}(Q_{j})\in\mathrm{Pth}_{\boldsymbol{\mathcal{A}}^{(2)},s_{j}}
\tag{2}
\\&\Leftrightarrow
\mathrm{ip}^{(2,X)@}_{s}(Q)\in\mathrm{Pth}_{\boldsymbol{\mathcal{A}}^{(2)},s}.
\tag{3}
\end{align*}

In the just stated chain of equivalences, the first equivalence follows from left to right by Lemma~\ref{LDThetaCongSub} and from right to left because $\sigma$ is a total operation in $\mathbf{Pth}_{\boldsymbol{\mathcal{A}}^{(2)}}$, according to Proposition~\ref{PDPthAlg}, and $\mathrm{ip}^{(2,X)@}$ is a $\Sigma^{\boldsymbol{\mathcal{A}}^{(2)}}$-homomorphism, according to Definition~\ref{DDIp}; the second equivalence follows by induction, finally the third equivalence follows from left to right because $\sigma$ is a total operation in $\mathbf{Pth}_{\boldsymbol{\mathcal{A}}^{(2)}}$, according to Proposition~\ref{PDPthAlg},  and $\mathrm{ip}^{(2,X)@}$ is a $\Sigma^{\boldsymbol{\mathcal{A}}^{(2)}}$-homomorphism, according to Definition~\ref{DDIp}, and from right to left by Lemma~\ref{LDThetaCongSub}.

Assume, without loss of generality, that $\mathrm{ip}^{(2,X)@}_{s}(P)$ is a second-order path in $\mathrm{Pth}_{\boldsymbol{\mathcal{A}}^{(2)},s}$. As we have seen before, this is the case exactly when, for every $j\in\bb{\mathbf{s}}$, $\mathrm{ip}^{(2,X)@}_{s_{j}}(P_{j})$ is a second-order path in $\mathrm{Pth}_{\boldsymbol{\mathcal{A}}^{(2)},s_{j}}$. By induction, we also have that, for every $j\in\bb{\mathbf{s}}$, $\mathrm{ip}^{(2,X)@}_{s_{j}}(Q_{j})$ is a second-order path in $\mathrm{Pth}_{\boldsymbol{\mathcal{A}}^{(2)},s_{j}}$. Moreover, by induction, the following equality holds
\[
\mathrm{CH}^{(2)}_{s_{j}}\left(
\mathrm{ip}^{(2,X)@}_{s_{j}}\left(
P_{j}
\right)\right)
=
\mathrm{CH}^{(2)}_{s_{j}}\left(
\mathrm{ip}^{(2,X)@}_{s_{j}}\left(
Q_{j}
\right)\right).
\]

Note that the following chain of equalities holds
\allowdisplaybreaks
\begin{align*}
\mathrm{CH}^{(2)}_{s}\left(
\mathrm{ip}^{(2,X)@}_{s}\left(
P
\right)
\right)
&=
\mathrm{CH}^{(2)}_{s}\left(
\mathrm{ip}^{(2,X)@}_{s}\left(
\sigma^{\mathbf{T}_{\Sigma^{\boldsymbol{\mathcal{A}}^{(2)}}}(X)}
\left(
\left(
P_{j}
\right)_{j\in\bb{\mathbf{s}}}
\right)
\right)
\right)
\tag{1}
\\&=
\mathrm{CH}^{(2)}_{s}\left(
\sigma^{\mathbf{Pth}_{\boldsymbol{\mathcal{A}}^{(2)}}}
\left(
\left(
\mathrm{ip}^{(2,X)@}_{s_{j}}\left(
P_{j}
\right)
\right)_{j\in\bb{\mathbf{s}}}
\right)
\right)
\tag{2}
\\&=
\mathrm{CH}^{(2)}_{s}\left(
\sigma^{\mathbf{Pth}_{\boldsymbol{\mathcal{A}}^{(2)}}}
\left(
\left(
\mathrm{ip}^{(2,X)@}_{s_{j}}\left(
Q_{j}
\right)
\right)_{j\in\bb{\mathbf{s}}}
\right)
\right)
\tag{3}
\\&=
\mathrm{CH}^{(2)}_{s}\left(
\mathrm{ip}^{(2,X)@}_{s}\left(
\sigma^{\mathbf{T}_{\Sigma^{\boldsymbol{\mathcal{A}}^{(2)}}}(X)}
\left(
\left(
Q_{j}
\right)_{j\in\bb{\mathbf{s}}}
\right)
\right)
\right)
\tag{4}
\\&=
\mathrm{CH}^{(2)}_{s}\left(
\mathrm{ip}^{(2,X)@}_{s}\left(
Q
\right)
\right).
\tag{5}
\end{align*}

In the just stated chain of equalities, the first equality unravels the description of $P$; the second equivalence follows from the fact that $\sigma$ is a total operation in $\mathbf{Pth}_{\boldsymbol{\mathcal{A}}^{(2)}}$, according to Proposition~\ref{PDPthAlg},  and $\mathrm{ip}^{(2,X)@}$ is a $\Sigma^{\boldsymbol{\mathcal{A}}^{(2)}}$-homomorphism, according to Definition~\ref{DDIp}, and the fact that we are assuming that $\mathrm{ip}^{(2,X)@}_{s}(P)$ is a second-order path in $\mathrm{Pth}_{\boldsymbol{\mathcal{A}}^{(2)},s}$; the third equality follows by induction and from the fact that, according to Proposition~\ref{PDCHCong}, $\mathrm{Ker}(\mathrm{CH}^{(2)})$ is a $\Sigma^{\boldsymbol{\mathcal{A}}^{(2)}}$-congruence; the fourth equality follows from the fact that $\sigma$ is a total operation in $\mathbf{Pth}_{\boldsymbol{\mathcal{A}}^{(2)}}$, according to Proposition~\ref{PDPthAlg},  and $\mathrm{ip}^{(2,X)@}$ is a $\Sigma^{\boldsymbol{\mathcal{A}}^{(2)}}$-homomorphism, according to Definition~\ref{DDIp}, and the fact that we are assuming that, for every $j\in\bb{\mathbf{s}}$, $\mathrm{ip}^{(2,X)@}_{s_{j}}(Q_{j})$ is a second-order path in $\mathrm{Pth}_{\boldsymbol{\mathcal{A}}^{(2)},s_{j}}$; finally, the last equality recovers the description of $Q$.

This completes the case $\sigma\in\Sigma_{\mathbf{s},s}$.

\textsf{$\gamma$ is the $0$-source operation symbol.}

Let $\gamma$ be the operation symbol $\mathrm{sc}^{0}_{s}$ in $\Sigma^{\boldsymbol{\mathcal{A}}}_{s,s}$ and let $(P',Q')$ be a family of pairs in $\mathrm{C}^{n}_{\Sigma^{\boldsymbol{\mathcal{A}}^{(2)}}}(\Theta^{(2)})_{s}$ for which 
$$
(P,Q)=
\left(
\mathrm{sc}^{0\mathbf{T}_{\Sigma^{\boldsymbol{\mathcal{A}}^{(2)}}}(X)}_{s}\left(
P'
\right)
,
\mathrm{sc}^{0\mathbf{T}_{\Sigma^{\boldsymbol{\mathcal{A}}^{(2)}}}(X)}_{s}\left(
Q'
\right)
\right).
$$

Note that the following chain of equivalences holds
\allowdisplaybreaks
\begin{align*}
\mathrm{ip}^{(2,X)@}_{s}(P)\in\mathrm{Pth}_{\boldsymbol{\mathcal{A}}^{(2)},s}
&\Leftrightarrow
\mathrm{ip}^{(2,X)@}_{s}(P')\in\mathrm{Pth}_{\boldsymbol{\mathcal{A}}^{(2)},s}
\tag{1}
\\&\Leftrightarrow
\mathrm{ip}^{(2,X)@}_{s}(Q')\in\mathrm{Pth}_{\boldsymbol{\mathcal{A}}^{(2)},s}
\tag{2}
\\&\Leftrightarrow
\mathrm{ip}^{(2,X)@}_{s}(Q)\in\mathrm{Pth}_{\boldsymbol{\mathcal{A}}^{(2)},s}.
\tag{3}
\end{align*}

In the just stated chain of equivalences, the first equivalence follows from left to right by Lemma~\ref{LDThetaCongSub} and from right to left because $\mathrm{sc}^{0}_{s}$ is a total operation in $\mathbf{Pth}_{\boldsymbol{\mathcal{A}}^{(2)}}$, according to Proposition~\ref{PDPthCatAlg}, and $\mathrm{ip}^{(2,X)@}$ is a $\Sigma^{\boldsymbol{\mathcal{A}}^{(2)}}$-homomorphism, according to Definition~\ref{DDIp}; the second equivalence follows by induction, finally the third equivalence follows from left to right because $\mathrm{sc}^{0}_{s}$ is a total operation in $\mathbf{Pth}_{\boldsymbol{\mathcal{A}}^{(2)}}$, according to Proposition~\ref{PDPthAlg},  and $\mathrm{ip}^{(2,X)@}$ is a $\Sigma^{\boldsymbol{\mathcal{A}}^{(2)}}$-homomorphism, according to Definition~\ref{DDIp}, and from right to left by Lemma~\ref{LDThetaCongSub}.

Assume, without loss of generality, that $\mathrm{ip}^{(2,X)@}_{s}(P)$ is a second-order path in $\mathrm{Pth}_{\boldsymbol{\mathcal{A}}^{(2)},s}$. As we have seen before, this is the case exactly when $\mathrm{ip}^{(2,X)@}_{s}(P')$ is a second-order path in $\mathrm{Pth}_{\boldsymbol{\mathcal{A}}^{(2)},s}$. By induction, we also have that $\mathrm{ip}^{(2,X)@}_{s}(Q')$ is a second-order path in $\mathrm{Pth}_{\boldsymbol{\mathcal{A}}^{(2)},s}$. Moreover, by induction, the following equality holds
\[
\mathrm{CH}^{(2)}_{s}\left(
\mathrm{ip}^{(2,X)@}_{s}\left(
P'
\right)\right)
=
\mathrm{CH}^{(2)}_{s}\left(
\mathrm{ip}^{(2,X)@}_{s}\left(
Q'
\right)\right).
\]

Note that the following chain of equalities holds
\allowdisplaybreaks
\begin{align*}
\mathrm{CH}^{(2)}_{s}\left(
\mathrm{ip}^{(2,X)@}_{s}\left(
P
\right)
\right)
&=
\mathrm{CH}^{(2)}_{s}\left(
\mathrm{ip}^{(2,X)@}_{s}\left(
\mathrm{sc}^{0\mathbf{T}_{\Sigma^{\boldsymbol{\mathcal{A}}^{(2)}}}(X)}_{s}
\left(
P'
\right)
\right)
\right)
\tag{1}
\\&=
\mathrm{CH}^{(2)}_{s}\left(
\mathrm{sc}^{0\mathbf{Pth}_{\boldsymbol{\mathcal{A}}^{(2)}}}_{s}
\left(
\mathrm{ip}^{(2,X)@}_{s}\left(
P'
\right)
\right)
\right)
\tag{2}
\\&=
\mathrm{CH}^{(2)}_{s}\left(
\mathrm{sc}^{0\mathbf{Pth}_{\boldsymbol{\mathcal{A}}^{(2)}}}_{s}
\left(
\mathrm{ip}^{(2,X)@}_{s}\left(
Q'
\right)
\right)
\right)
\tag{3}
\\&=
\mathrm{CH}^{(2)}_{s}\left(
\mathrm{ip}^{(2,X)@}_{s}\left(
\mathrm{sc}^{0\mathbf{T}_{\Sigma^{\boldsymbol{\mathcal{A}}^{(2)}}}(X)}_{s}
\left(
Q'
\right)
\right)
\right)
\tag{4}
\\&=
\mathrm{CH}^{(2)}_{s}\left(
\mathrm{ip}^{(2,X)@}_{s}\left(
Q
\right)
\right).
\tag{5}
\end{align*}

In the just stated chain of equalities, the first equality unravels the description of $P$; the second equivalence follows from the fact that $\mathrm{sc}^{0}_{s}$ is a total operation in $\mathbf{Pth}_{\boldsymbol{\mathcal{A}}^{(2)}}$, according to Proposition~\ref{PDPthCatAlg},  and $\mathrm{ip}^{(2,X)@}$ is a $\Sigma^{\boldsymbol{\mathcal{A}}^{(2)}}$-homomorphism, according to Definition~\ref{DDIp}, and the fact that we are assuming that $\mathrm{ip}^{(2,X)@}_{s}(P)$ is a second-order path in $\mathrm{Pth}_{\boldsymbol{\mathcal{A}}^{(2)},s}$; the third equality follows by induction and from the fact that, according to Proposition~\ref{PDCHCong}, $\mathrm{Ker}(\mathrm{CH}^{(2)})$ is a $\Sigma^{\boldsymbol{\mathcal{A}}^{(2)}}$-congruence; the fourth equality follows from the fact that $\mathrm{sc}^{0}_{s}$ is a total operation in $\mathbf{Pth}_{\boldsymbol{\mathcal{A}}^{(2)}}$, according to Proposition~\ref{PDPthCatAlg},  and $\mathrm{ip}^{(2,X)@}$ is a $\Sigma^{\boldsymbol{\mathcal{A}}^{(2)}}$-homomorphism, according to Definition~\ref{DDIp}, and the fact that we are assuming that $\mathrm{ip}^{(2,X)@}_{s}(Q')$ is a second-order path in $\mathrm{Pth}_{\boldsymbol{\mathcal{A}}^{(2)},s}$; finally, the last equality recovers the description of $Q$.

This completes the case of the $0$-source.

\textsf{$\gamma$ is the $0$-target operation symbol.}

Let $\gamma$ be the operation symbol $\mathrm{tg}^{0}_{s}$ in $\Sigma^{\boldsymbol{\mathcal{A}}}_{s,s}$ and let $(P',Q')$ be a family of pairs in $\mathrm{C}^{n}_{\Sigma^{\boldsymbol{\mathcal{A}}^{(2)}}}(\Theta^{(2)})_{s}$ for which 
$$
(P,Q)=
\left(
\mathrm{tg}^{0\mathbf{T}_{\Sigma^{\boldsymbol{\mathcal{A}}^{(2)}}}(X)}_{s}\left(
P'
\right)
,
\mathrm{tg}^{0\mathbf{T}_{\Sigma^{\boldsymbol{\mathcal{A}}^{(2)}}}(X)}_{s}\left(
Q'
\right)
\right).
$$

Then, the following properties hold
\begin{itemize}
\item[(i)] $\mathrm{ip}^{(2,X)@}_{s}(P)\in\mathrm{Pth}_{\boldsymbol{\mathcal{A}}^{(2)},s}$ if, and only if, $\mathrm{ip}^{(2,X)@}_{s}(Q)\in\mathrm{Pth}_{\boldsymbol{\mathcal{A}}^{(2)},s}$;
\item[(ii)] If $\mathrm{ip}^{(2,X)@}_{s}(P)\in\mathrm{Pth}_{\boldsymbol{\mathcal{A}}^{(2)},s}$ or  $\mathrm{ip}^{(2,X)@}_{s}(Q)\in\mathrm{Pth}_{\boldsymbol{\mathcal{A}}^{(2)},s}$ is a second-order path in $\mathrm{Pth}_{\boldsymbol{\mathcal{A}}^{(2)},s}$ then
$
\mathrm{CH}^{(2)}_{s}(
\mathrm{ip}^{(2,X)@}_{s}(
P))=\mathrm{CH}^{(2)}_{s}(
\mathrm{ip}^{(2,X)@}_{s}(
Q)).
$
\end{itemize}

The proof of this case is similar to that of the $0$-source.

This completes the case of the $0$-target.

\textsf{$\gamma$ is the $0$-composition operation symbol.}

Let $\gamma$ be the operation symbol $\circ^{0}_{s}$ in $\Sigma^{\boldsymbol{\mathcal{A}}}_{s,s}$ and let $(P',Q')$ and $(P'',Q'')$ be two families of pairs in $\mathrm{C}^{n}_{\Sigma^{\boldsymbol{\mathcal{A}}^{(2)}}}(\Theta^{(2)})_{s}$ for which 
$$
(P,Q)=
\left(
P''
\circ^{0\mathbf{T}_{\Sigma^{\boldsymbol{\mathcal{A}}^{(2)}}}(X)}_{s}
P'
,
Q''
\circ^{0\mathbf{T}_{\Sigma^{\boldsymbol{\mathcal{A}}^{(2)}}}(X)}_{s}
Q'
\right).
$$

Note that the following chain of equivalences holds
\allowdisplaybreaks
\begin{align*}
\mathrm{ip}^{(2,X)@}_{s}(P)\in\mathrm{Pth}_{\boldsymbol{\mathcal{A}}^{(2)},s}
&\Leftrightarrow
\left\lbrace
\begin{array}{l}
\mathrm{ip}^{(2,X)@}_{s}(P'),\mathrm{ip}^{(2,X)@}_{s}(P'')\in\mathrm{Pth}_{\boldsymbol{\mathcal{A}}^{(2)},s}
\\
\mathrm{sc}^{(0,2)}_{s}(
\mathrm{ip}^{(2,X)@}_{s}(P'')
)
=
\mathrm{tg}^{(0,2)}_{s}(
\mathrm{ip}^{(2,X)@}_{s}(P')
)
\end{array}
\right.
\tag{1}
\\&\Leftrightarrow
\left\lbrace
\begin{array}{l}
\mathrm{ip}^{(2,X)@}_{s}(Q'),\mathrm{ip}^{(2,X)@}_{s}(Q'')\in\mathrm{Pth}_{\boldsymbol{\mathcal{A}}^{(2)},s}
\\
\mathrm{sc}^{(0,2)}_{s}(
\mathrm{ip}^{(2,X)@}_{s}(Q'')
)
=
\mathrm{tg}^{(0,2)}_{s}(
\mathrm{ip}^{(2,X)@}_{s}(Q')
)
\end{array}
\right.
\tag{2}
\\&\Leftrightarrow
\mathrm{ip}^{(2,X)@}_{s}(Q)\in\mathrm{Pth}_{\boldsymbol{\mathcal{A}}^{(2)},s}.
\tag{3}
\end{align*}

In the just stated chain of equivalences, the first equivalence follows from left to right by Lemma~\ref{LDThetaCongSub} and by the description of the $0$-composition operation in $\mathbf{Pth}_{\boldsymbol{\mathcal{A}}^{(2)}}$, according to Proposition~\ref{PDPthCatAlg}, and from right to left because $\circ^{0}_{s}$ is a defined operation in $\mathbf{Pth}_{\boldsymbol{\mathcal{A}}^{(2)}}$, according to Proposition~\ref{PDPthCatAlg}, and $\mathrm{ip}^{(2,X)@}$ is a $\Sigma^{\boldsymbol{\mathcal{A}}^{(2)}}$-homomorphism, according to Definition~\ref{DDIp}; the second equivalence follows by induction. Note also that
\allowdisplaybreaks
\begin{align*}
\mathrm{CH}^{(2)}_{s}\left(
\mathrm{ip}^{(2,X)@}_{s}\left(
P'
\right)
\right)
&=
\mathrm{CH}^{(2)}_{s}\left(
\mathrm{ip}^{(2,X)@}_{s}\left(
Q'
\right)
\right);
\\
\mathrm{CH}^{(2)}_{s}\left(
\mathrm{ip}^{(2,X)@}_{s}\left(
P''
\right)
\right)
&=
\mathrm{CH}^{(2)}_{s}\left(
\mathrm{ip}^{(2,X)@}_{s}\left(
Q''
\right)
\right).
\end{align*}
Therefore, according to Corollary~\ref{CDCH}, we have that 
\allowdisplaybreaks
\begin{align*}
\mathrm{sc}^{(0,2)}_{s}\left(
\mathrm{ip}^{(2,X)@}_{s}\left(
P''
\right)\right)&=
\mathrm{sc}^{(0,2)}_{s}\left(
\mathrm{ip}^{(2,X)@}_{s}\left(
Q''
\right)\right)
\\
\mathrm{tg}^{(0,2)}_{s}\left(
\mathrm{ip}^{(2,X)@}_{s}\left(
P'
\right)\right)&=
\mathrm{tg}^{(0,2)}_{s}\left(
\mathrm{ip}^{(2,X)@}_{s}\left(
Q'
\right)\right);
\end{align*}
finally the third equivalence follows from left to right because $\circ^{0}_{s}$ is a defined operation in $\mathbf{Pth}_{\boldsymbol{\mathcal{A}}^{(2)}}$, according to Proposition~\ref{PDPthCatAlg}, and $\mathrm{ip}^{(2,X)@}$ is a $\Sigma^{\boldsymbol{\mathcal{A}}^{(2)}}$-homomorphism, according to Definition~\ref{DDIp}, and from right to left by Lemma~\ref{LDThetaCongSub} and by the description of the $0$-composition operation in $\mathbf{Pth}_{\boldsymbol{\mathcal{A}}^{(2)}}$, according to Proposition~\ref{PDPthCatAlg}.

Assume, without loss of generality, that $\mathrm{ip}^{(2,X)@}_{s}(P)$ is a second-order path in $\mathrm{Pth}_{\boldsymbol{\mathcal{A}}^{(2)},s}$. As we have seen before, this is the case exactly when $\mathrm{ip}^{(2,X)@}_{s}(P')$ and $\mathrm{ip}^{(2,X)@}_{s}(P'')$ are second-order paths in $\mathrm{Pth}_{\boldsymbol{\mathcal{A}}^{(2)},s}$ satisfying that 
\[
\mathrm{sc}^{(0,2)}_{s}\left(
\mathrm{ip}^{(2,X)@}_{s}\left(
P''
\right)\right)=
\mathrm{tg}^{(0,2)}_{s}\left(
\mathrm{ip}^{(2,X)@}_{s}\left(
P'
\right)\right).
\]
By induction, we also have that $\mathrm{ip}^{(2,X)@}_{s}(Q')$ and $\mathrm{ip}^{(2,X)@}_{s}(Q'')$ are second-order paths in $\mathrm{Pth}_{\boldsymbol{\mathcal{A}}^{(2)},s}$. As we have proven before, their $0$-composition is well-defined. Moreover, by induction, the following equality holds
\allowdisplaybreaks
\begin{align*}
\mathrm{CH}^{(2)}_{s}\left(
\mathrm{ip}^{(2,X)@}_{s}\left(
P'
\right)\right)
&=
\mathrm{CH}^{(2)}_{s}\left(
\mathrm{ip}^{(2,X)@}_{s}\left(
Q'
\right)\right);
\\
\mathrm{CH}^{(2)}_{s}\left(
\mathrm{ip}^{(2,X)@}_{s}\left(
P''
\right)\right)
&=
\mathrm{CH}^{(2)}_{s}\left(
\mathrm{ip}^{(2,X)@}_{s}\left(
Q''
\right)\right);
\end{align*}

Note that the following chain of equalities holds
\allowdisplaybreaks
\begin{align*}
\mathrm{CH}^{(2)}_{s}\left(
\mathrm{ip}^{(2,X)@}_{s}\left(
P
\right)
\right)
&=
\mathrm{CH}^{(2)}_{s}\left(
\mathrm{ip}^{(2,X)@}_{s}\left(
P''
\circ^{0\mathbf{T}_{\Sigma^{\boldsymbol{\mathcal{A}}^{(2)}}}(X)}_{s}
P'
\right)
\right)
\tag{1}
\\&=
\mathrm{CH}^{(2)}_{s}\left(
\mathrm{ip}^{(2,X)@}_{s}\left(
P''
\right)
\circ^{0\mathbf{Pth}_{\boldsymbol{\mathcal{A}}^{(2)}}}_{s}
\mathrm{ip}^{(2,X)@}_{s}\left(
P'
\right)
\right)
\tag{2}
\\&=
\mathrm{CH}^{(2)}_{s}\left(
\mathrm{ip}^{(2,X)@}_{s}\left(
Q''
\right)
\circ^{0\mathbf{Pth}_{\boldsymbol{\mathcal{A}}^{(2)}}}_{s}
\mathrm{ip}^{(2,X)@}_{s}\left(
Q'
\right)
\right)
\tag{3}
\\&=
\mathrm{CH}^{(2)}_{s}\left(
\mathrm{ip}^{(2,X)@}_{s}\left(
Q''
\circ^{0\mathbf{T}_{\Sigma^{\boldsymbol{\mathcal{A}}^{(2)}}}(X)}_{s}
Q'
\right)
\right)
\tag{4}
\\&=
\mathrm{CH}^{(2)}_{s}\left(
\mathrm{ip}^{(2,X)@}_{s}\left(
Q
\right)
\right).
\tag{5}
\end{align*}

In the just stated chain of equalities, the first equality unravels the description of $P$; the second equivalence follows from the fact that $\circ^{0}_{s}$ is a defined operation in $\mathbf{Pth}_{\boldsymbol{\mathcal{A}}^{(2)}}$, according to Proposition~\ref{PDPthCatAlg},  and $\mathrm{ip}^{(2,X)@}$ is a $\Sigma^{\boldsymbol{\mathcal{A}}^{(2)}}$-homomorphism, according to Definition~\ref{DDIp}, and the fact that we are assuming that $\mathrm{ip}^{(2,X)@}_{s}(P)$ is a second-order path in $\mathrm{Pth}_{\boldsymbol{\mathcal{A}}^{(2)},s}$; the third equality follows by induction and from the fact that, according to Proposition~\ref{PDCHCong}, $\mathrm{Ker}(\mathrm{CH}^{(2)})$ is a $\Sigma^{\boldsymbol{\mathcal{A}}^{(2)}}$-congruence; the fourth equality follows from the fact that $\circ^{0}_{s}$ is a defined operation in $\mathbf{Pth}_{\boldsymbol{\mathcal{A}}^{(2)}}$, according to Proposition~\ref{PDPthCatAlg},  and $\mathrm{ip}^{(2,X)@}$ is a $\Sigma^{\boldsymbol{\mathcal{A}}^{(2)}}$-homomorphism, according to Definition~\ref{DDIp}, and the fact that we are assuming that $\mathrm{ip}^{(2,X)@}_{s}(Q')$ and $\mathrm{ip}^{(2,X)@}_{s}(Q'')$ are second-order paths in $\mathrm{Pth}_{\boldsymbol{\mathcal{A}}^{(2)},s}$ whose $0$-composition is well-defined; finally, the last equality recovers the description of $Q$.

This completes the case of the $0$-composition.

\textsf{$\gamma$ is the $1$-source operation symbol.}

Let $\gamma$ be the operation symbol $\mathrm{sc}^{1}_{s}$ in $\Sigma^{\boldsymbol{\mathcal{A}}^{(2)}}_{s,s}$ and let $(P',Q')$ be a family of pairs in $\mathrm{C}^{n}_{\Sigma^{\boldsymbol{\mathcal{A}}^{(2)}}}(\Theta^{(2)})_{s}$ for which 
$$
(P,Q)=
\left(
\mathrm{sc}^{1\mathbf{T}_{\Sigma^{\boldsymbol{\mathcal{A}}^{(2)}}}(X)}_{s}\left(
P'
\right)
,
\mathrm{sc}^{1\mathbf{T}_{\Sigma^{\boldsymbol{\mathcal{A}}^{(2)}}}(X)}_{s}\left(
Q'
\right)
\right).
$$

Note that the following chain of equivalences holds
\allowdisplaybreaks
\begin{align*}
\mathrm{ip}^{(2,X)@}_{s}(P)\in\mathrm{Pth}_{\boldsymbol{\mathcal{A}}^{(2)},s}
&\Leftrightarrow
\mathrm{ip}^{(2,X)@}_{s}(P')\in\mathrm{Pth}_{\boldsymbol{\mathcal{A}}^{(2)},s}
\tag{1}
\\&\Leftrightarrow
\mathrm{ip}^{(2,X)@}_{s}(Q')\in\mathrm{Pth}_{\boldsymbol{\mathcal{A}}^{(2)},s}
\tag{2}
\\&\Leftrightarrow
\mathrm{ip}^{(2,X)@}_{s}(Q)\in\mathrm{Pth}_{\boldsymbol{\mathcal{A}}^{(2)},s}.
\tag{3}
\end{align*}

In the just stated chain of equivalences, the first equivalence follows from left to right by Lemma~\ref{LDThetaCongSub} and from right to left because $\mathrm{sc}^{1}_{s}$ is a total operation in $\mathbf{Pth}_{\boldsymbol{\mathcal{A}}^{(2)}}$, according to Proposition~\ref{PDPthDCatAlg}, and $\mathrm{ip}^{(2,X)@}$ is a $\Sigma^{\boldsymbol{\mathcal{A}}^{(2)}}$-homomorphism, according to Definition~\ref{DDIp}; the second equivalence follows by induction, finally the third equivalence follows from left to right because $\mathrm{sc}^{1}_{s}$ is a total operation in $\mathbf{Pth}_{\boldsymbol{\mathcal{A}}^{(2)}}$, according to Proposition~\ref{PDPthAlg},  and $\mathrm{ip}^{(2,X)@}$ is a $\Sigma^{\boldsymbol{\mathcal{A}}^{(2)}}$-homomorphism, according to Definition~\ref{DDIp}, and from right to left by Lemma~\ref{LDThetaCongSub}.

Assume, without loss of generality, that $\mathrm{ip}^{(2,X)@}_{s}(P)$ is a second-order path in $\mathrm{Pth}_{\boldsymbol{\mathcal{A}}^{(2)},s}$. As we have seen before, this is the case exactly when $\mathrm{ip}^{(2,X)@}_{s}(P')$ is a second-order path in $\mathrm{Pth}_{\boldsymbol{\mathcal{A}}^{(2)},s}$. By induction, we also have that $\mathrm{ip}^{(2,X)@}_{s}(Q')$ is a second-order path in $\mathrm{Pth}_{\boldsymbol{\mathcal{A}}^{(2)},s}$. Moreover, by induction, the following equality holds
\[
\mathrm{CH}^{(2)}_{s}\left(
\mathrm{ip}^{(2,X)@}_{s}\left(
P'
\right)\right)
=
\mathrm{CH}^{(2)}_{s}\left(
\mathrm{ip}^{(2,X)@}_{s}\left(
Q'
\right)\right).
\]

Note that the following chain of equalities holds
\allowdisplaybreaks
\begin{align*}
\mathrm{CH}^{(2)}_{s}\left(
\mathrm{ip}^{(2,X)@}_{s}\left(
P
\right)
\right)
&=
\mathrm{CH}^{(2)}_{s}\left(
\mathrm{ip}^{(2,X)@}_{s}\left(
\mathrm{sc}^{1\mathbf{T}_{\Sigma^{\boldsymbol{\mathcal{A}}^{(2)}}}(X)}_{s}
\left(
P'
\right)
\right)
\right)
\tag{1}
\\&=
\mathrm{CH}^{(2)}_{s}\left(
\mathrm{sc}^{1\mathbf{Pth}_{\boldsymbol{\mathcal{A}}^{(2)}}}_{s}
\left(
\mathrm{ip}^{(2,X)@}_{s}\left(
P'
\right)
\right)
\right)
\tag{2}
\\&=
\mathrm{CH}^{(2)}_{s}\left(
\mathrm{sc}^{1\mathbf{Pth}_{\boldsymbol{\mathcal{A}}^{(2)}}}_{s}
\left(
\mathrm{ip}^{(2,X)@}_{s}\left(
Q'
\right)
\right)
\right)
\tag{3}
\\&=
\mathrm{CH}^{(2)}_{s}\left(
\mathrm{ip}^{(2,X)@}_{s}\left(
\mathrm{sc}^{1\mathbf{T}_{\Sigma^{\boldsymbol{\mathcal{A}}^{(2)}}}(X)}_{s}
\left(
Q'
\right)
\right)
\right)
\tag{4}
\\&=
\mathrm{CH}^{(2)}_{s}\left(
\mathrm{ip}^{(2,X)@}_{s}\left(
Q
\right)
\right).
\tag{5}
\end{align*}

In the just stated chain of equalities, the first equality unravels the description of $P$; the second equivalence follows from the fact that $\mathrm{sc}^{1}_{s}$ is a total operation in $\mathbf{Pth}_{\boldsymbol{\mathcal{A}}^{(2)}}$, according to Proposition~\ref{PDPthDCatAlg},  and $\mathrm{ip}^{(2,X)@}$ is a $\Sigma^{\boldsymbol{\mathcal{A}}^{(2)}}$-homomorphism, according to Definition~\ref{DDIp}, and the fact that we are assuming that $\mathrm{ip}^{(2,X)@}_{s}(P)$ is a second-order path in $\mathrm{Pth}_{\boldsymbol{\mathcal{A}}^{(2)},s}$; the third equality follows by induction and from the fact that, according to Proposition~\ref{PDCHCong}, $\mathrm{Ker}(\mathrm{CH}^{(2)})$ is a $\Sigma^{\boldsymbol{\mathcal{A}}^{(2)}}$-congruence; the fourth equality follows from the fact that $\mathrm{sc}^{1}_{s}$ is a total operation in $\mathbf{Pth}_{\boldsymbol{\mathcal{A}}^{(2)}}$, according to Proposition~\ref{PDPthDCatAlg},  and $\mathrm{ip}^{(2,X)@}$ is a $\Sigma^{\boldsymbol{\mathcal{A}}^{(2)}}$-homomorphism, according to Definition~\ref{DDIp}, and the fact that we are assuming that $\mathrm{ip}^{(2,X)@}_{s}(Q')$ is a second-order path in $\mathrm{Pth}_{\boldsymbol{\mathcal{A}}^{(2)},s}$; finally, the last equality recovers the description of $Q$.

This completes the case of the $1$-source.

\textsf{$\gamma$ is the $1$-target operation symbol.}

Let $\gamma$ be the operation symbol $\mathrm{tg}^{1}_{s}$ in $\Sigma^{\boldsymbol{\mathcal{A}}^{(2)}}_{s,s}$ and let $(P',Q')$ be a family of pairs in $\mathrm{C}^{n}_{\Sigma^{\boldsymbol{\mathcal{A}}^{(2)}}}(\Theta^{(2)})_{s}$ for which 
$$
(P,Q)=
\left(
\mathrm{tg}^{1\mathbf{T}_{\Sigma^{\boldsymbol{\mathcal{A}}^{(2)}}}(X)}_{s}\left(
P'
\right)
,
\mathrm{tg}^{1\mathbf{T}_{\Sigma^{\boldsymbol{\mathcal{A}}^{(2)}}}(X)}_{s}\left(
Q'
\right)
\right).
$$

Then, the following properties hold
\begin{itemize}
\item[(i)] $\mathrm{ip}^{(2,X)@}_{s}(P)\in\mathrm{Pth}_{\boldsymbol{\mathcal{A}}^{(2)},s}$ if, and only if, $\mathrm{ip}^{(2,X)@}_{s}(Q)\in\mathrm{Pth}_{\boldsymbol{\mathcal{A}}^{(2)},s}$;
\item[(ii)] If $\mathrm{ip}^{(2,X)@}_{s}(P)\in\mathrm{Pth}_{\boldsymbol{\mathcal{A}}^{(2)},s}$ or  $\mathrm{ip}^{(2,X)@}_{s}(Q)\in\mathrm{Pth}_{\boldsymbol{\mathcal{A}}^{(2)},s}$ is a second-order path in $\mathrm{Pth}_{\boldsymbol{\mathcal{A}}^{(2)},s}$ then
$
\mathrm{CH}^{(2)}_{s}(
\mathrm{ip}^{(2,X)@}_{s}(
P))=\mathrm{CH}^{(2)}_{s}(
\mathrm{ip}^{(2,X)@}_{s}(
Q)).
$
\end{itemize}

The proof of this case is similar to that of the $1$-source.

This completes the case of the $1$-target.

\textsf{$\gamma$ is the $1$-composition operation symbol.}

Let $\gamma$ be the operation symbol $\circ^{1}_{s}$ in $\Sigma^{\boldsymbol{\mathcal{A}}^{(2)}}_{s,s}$ and let $(P',Q')$ and $(P'',Q'')$ be two families of pairs in $\mathrm{C}^{n}_{\Sigma^{\boldsymbol{\mathcal{A}}^{(2)}}}(\Theta^{(2)})_{s}$ for which 
$$
(P,Q)=
\left(
P''
\circ^{1\mathbf{T}_{\Sigma^{\boldsymbol{\mathcal{A}}^{(2)}}}(X)}_{s}
P'
,
Q''
\circ^{1\mathbf{T}_{\Sigma^{\boldsymbol{\mathcal{A}}^{(2)}}}(X)}_{s}
Q'
\right).
$$

Note that the following chain of equivalences holds
\allowdisplaybreaks
\begin{align*}
\mathrm{ip}^{(2,X)@}_{s}(P)\in\mathrm{Pth}_{\boldsymbol{\mathcal{A}}^{(2)},s}
&\Leftrightarrow
\left\lbrace
\begin{array}{l}
\mathrm{ip}^{(2,X)@}_{s}(P'),\mathrm{ip}^{(2,X)@}_{s}(P'')\in\mathrm{Pth}_{\boldsymbol{\mathcal{A}}^{(2)},s}
\\
\mathrm{sc}^{([1],2)}_{s}(
\mathrm{ip}^{(2,X)@}_{s}(P'')
)
=
\mathrm{tg}^{([1],2)}_{s}(
\mathrm{ip}^{(2,X)@}_{s}(P')
)
\end{array}
\right.
\tag{1}
\\&\Leftrightarrow
\left\lbrace
\begin{array}{l}
\mathrm{ip}^{(2,X)@}_{s}(Q'),\mathrm{ip}^{(2,X)@}_{s}(Q'')\in\mathrm{Pth}_{\boldsymbol{\mathcal{A}}^{(2)},s}
\\
\mathrm{sc}^{([1],2)}_{s}(
\mathrm{ip}^{(2,X)@}_{s}(Q'')
)
=
\mathrm{tg}^{([1],2)}_{s}(
\mathrm{ip}^{(2,X)@}_{s}(Q')
)
\end{array}
\right.
\tag{2}
\\&\Leftrightarrow
\mathrm{ip}^{(2,X)@}_{s}(Q)\in\mathrm{Pth}_{\boldsymbol{\mathcal{A}}^{(2)},s}.
\tag{3}
\end{align*}

In the just stated chain of equivalences, the first equivalence follows from left to right by Lemma~\ref{LDThetaCongSub} and by the description of the $1$-composition operation in $\mathbf{Pth}_{\boldsymbol{\mathcal{A}}^{(2)}}$, according to Proposition~\ref{PDPthDCatAlg}, and from right to left because $\circ^{1}_{s}$ is a defined operation in $\mathbf{Pth}_{\boldsymbol{\mathcal{A}}^{(2)}}$, according to Proposition~\ref{PDPthDCatAlg}, and $\mathrm{ip}^{(2,X)@}$ is a $\Sigma^{\boldsymbol{\mathcal{A}}^{(2)}}$-homomorphism, according to Definition~\ref{DDIp}; the second equivalence follows by induction. Note also that
\allowdisplaybreaks
\begin{align*}
\mathrm{CH}^{(2)}_{s}\left(
\mathrm{ip}^{(2,X)@}_{s}\left(
P'
\right)
\right)
&=
\mathrm{CH}^{(2)}_{s}\left(
\mathrm{ip}^{(2,X)@}_{s}\left(
Q'
\right)
\right);
\\
\mathrm{CH}^{(2)}_{s}\left(
\mathrm{ip}^{(2,X)@}_{s}\left(
P''
\right)
\right)
&=
\mathrm{CH}^{(2)}_{s}\left(
\mathrm{ip}^{(2,X)@}_{s}\left(
Q''
\right)
\right).
\end{align*}
Therefore, according to Lemma~\ref{LDCH}, we have that 
\allowdisplaybreaks
\begin{align*}
\mathrm{sc}^{([1],2)}_{s}\left(
\mathrm{ip}^{(2,X)@}_{s}\left(
P''
\right)\right)&=
\mathrm{sc}^{([1],2)}_{s}\left(
\mathrm{ip}^{(2,X)@}_{s}\left(
Q''
\right)\right)
\\
\mathrm{tg}^{([1],2)}_{s}\left(
\mathrm{ip}^{(2,X)@}_{s}\left(
P'
\right)\right)&=
\mathrm{tg}^{([1],2)}_{s}\left(
\mathrm{ip}^{(2,X)@}_{s}\left(
Q'
\right)\right);
\end{align*}
finally the third equivalence follows from left to right because $\circ^{1}_{s}$ is a defined operation in $\mathbf{Pth}_{\boldsymbol{\mathcal{A}}^{(2)}}$, according to Proposition~\ref{PDPthDCatAlg}, and $\mathrm{ip}^{(2,X)@}$ is a $\Sigma^{\boldsymbol{\mathcal{A}}^{(2)}}$-homomorphism, according to Definition~\ref{DDIp}, and from right to left by Lemma~\ref{LDThetaCongSub} and by the description of the $1$-composition operation in $\mathbf{Pth}_{\boldsymbol{\mathcal{A}}^{(2)}}$, according to Proposition~\ref{PDPthDCatAlg}.

Assume, without loss of generality, that $\mathrm{ip}^{(2,X)@}_{s}(P)$ is a second-order path in $\mathrm{Pth}_{\boldsymbol{\mathcal{A}}^{(2)},s}$. As we have seen before, this is the case exactly when $\mathrm{ip}^{(2,X)@}_{s}(P')$ and $\mathrm{ip}^{(2,X)@}_{s}(P'')$ are second-order paths in $\mathrm{Pth}_{\boldsymbol{\mathcal{A}}^{(2)},s}$ satisfying that 
\[
\mathrm{sc}^{([1],2)}_{s}\left(
\mathrm{ip}^{(2,X)@}_{s}\left(
P''
\right)\right)=
\mathrm{tg}^{([1],2)}_{s}\left(
\mathrm{ip}^{(2,X)@}_{s}\left(
P'
\right)\right).
\]
By induction, we also have that $\mathrm{ip}^{(2,X)@}_{s}(Q')$ and $\mathrm{ip}^{(2,X)@}_{s}(Q'')$ are second-order paths in $\mathrm{Pth}_{\boldsymbol{\mathcal{A}}^{(2)},s}$. As we have proven before, their $1$-composition is well-defined. Moreover, by induction, the following equality holds
\allowdisplaybreaks
\begin{align*}
\mathrm{CH}^{(2)}_{s}\left(
\mathrm{ip}^{(2,X)@}_{s}\left(
P'
\right)\right)
&=
\mathrm{CH}^{(2)}_{s}\left(
\mathrm{ip}^{(2,X)@}_{s}\left(
Q'
\right)\right);
\\
\mathrm{CH}^{(2)}_{s}\left(
\mathrm{ip}^{(2,X)@}_{s}\left(
P''
\right)\right)
&=
\mathrm{CH}^{(2)}_{s}\left(
\mathrm{ip}^{(2,X)@}_{s}\left(
Q''
\right)\right);
\end{align*}

Note that the following chain of equalities holds
\allowdisplaybreaks
\begin{align*}
\mathrm{CH}^{(2)}_{s}\left(
\mathrm{ip}^{(2,X)@}_{s}\left(
P
\right)
\right)
&=
\mathrm{CH}^{(2)}_{s}\left(
\mathrm{ip}^{(2,X)@}_{s}\left(
P''
\circ^{1\mathbf{T}_{\Sigma^{\boldsymbol{\mathcal{A}}^{(2)}}}(X)}_{s}
P'
\right)
\right)
\tag{1}
\\&=
\mathrm{CH}^{(2)}_{s}\left(
\mathrm{ip}^{(2,X)@}_{s}\left(
P''
\right)
\circ^{1\mathbf{Pth}_{\boldsymbol{\mathcal{A}}^{(2)}}}_{s}
\mathrm{ip}^{(2,X)@}_{s}\left(
P'
\right)
\right)
\tag{2}
\\&=
\mathrm{CH}^{(2)}_{s}\left(
\mathrm{ip}^{(2,X)@}_{s}\left(
Q''
\right)
\circ^{1\mathbf{Pth}_{\boldsymbol{\mathcal{A}}^{(2)}}}_{s}
\mathrm{ip}^{(2,X)@}_{s}\left(
Q'
\right)
\right)
\tag{3}
\\&=
\mathrm{CH}^{(2)}_{s}\left(
\mathrm{ip}^{(2,X)@}_{s}\left(
Q''
\circ^{1\mathbf{T}_{\Sigma^{\boldsymbol{\mathcal{A}}^{(2)}}}(X)}_{s}
Q'
\right)
\right)
\tag{4}
\\&=
\mathrm{CH}^{(2)}_{s}\left(
\mathrm{ip}^{(2,X)@}_{s}\left(
Q
\right)
\right).
\tag{5}
\end{align*}

In the just stated chain of equalities, the first equality unravels the description of $P$; the second equivalence follows from the fact that $\circ^{1}_{s}$ is a defined operation in $\mathbf{Pth}_{\boldsymbol{\mathcal{A}}^{(2)}}$, according to Proposition~\ref{PDPthDCatAlg},  and $\mathrm{ip}^{(2,X)@}$ is a $\Sigma^{\boldsymbol{\mathcal{A}}^{(2)}}$-homomorphism, according to Definition~\ref{DDIp}, and the fact that we are assuming that $\mathrm{ip}^{(2,X)@}_{s}(P)$ is a second-order path in $\mathrm{Pth}_{\boldsymbol{\mathcal{A}}^{(2)},s}$; the third equality follows by induction and from the fact that, according to Proposition~\ref{PDCHCong}, $\mathrm{Ker}(\mathrm{CH}^{(2)})$ is a $\Sigma^{\boldsymbol{\mathcal{A}}^{(2)}}$-congruence; the fourth equality follows from the fact that $\circ^{1}_{s}$ is a defined operation in $\mathbf{Pth}_{\boldsymbol{\mathcal{A}}^{(2)}}$, according to Proposition~\ref{PDPthCatAlg},  and $\mathrm{ip}^{(2,X)@}$ is a $\Sigma^{\boldsymbol{\mathcal{A}}^{(2)}}$-homomorphism, according to Definition~\ref{DDIp}, and the fact that we are assuming that $\mathrm{ip}^{(2,X)@}_{s}(Q')$ and $\mathrm{ip}^{(2,X)@}_{s}(Q'')$ are second-order paths in $\mathrm{Pth}_{\boldsymbol{\mathcal{A}}^{(2)},s}$ whose $1$-composition is well-defined; finally, the last equality recovers the description of $Q$.

This completes the case of the $1$-composition.

This completes the proof.
\end{proof}

\begin{restatable}{corollary}{CDThetaCong}
\label{CDThetaCong} Let $s$ be a sort in $S$, $\mathfrak{P}^{(2)}$ a second-order path in $\mathrm{Pth}_{\boldsymbol{\mathcal{A}}^{(2)},s}$ and $P$ a term in $\mathrm{T}_{\Sigma^{\boldsymbol{\mathcal{A}}^{(2)}}}(X)_{s}$ such that $(P,\mathrm{CH}^{(2)}(\mathfrak{P}^{(2)}))\in\Theta^{[2]}_{s}$. Then $\mathrm{ip}^{(2,X)@}_{s}(P)$ is a second-order path in $[\mathfrak{P}^{(2)}]^{}_{s}$. 
\end{restatable}
\begin{proof}
By Proposition~\ref{PDIpDCH}, $\mathrm{ip}^{(2,X)@}_{s}(\mathrm{CH}^{(2)}_{s}(\mathfrak{P}^{(2)}))$ is a second-order path in $[\mathfrak{P}^{(2)}]^{}_{s}$. Following Lemma~\ref{LDThetaCong}, we have that $\mathrm{ip}^{(2,X)@}_{s}(P)$ is a second-order path in $\mathrm{Pth}_{\boldsymbol{\mathcal{A}}^{(2)},s}$. Moreover, the following chain of equalities holds
$$
\mathrm{CH}^{(2)}_{s}\left(
\mathrm{ip}^{(2,X)@}_{s}\left(
P
\right)\right)
=
\mathrm{CH}^{(2)}_{s}\left(
\mathrm{ip}^{(2,X)@}_{s}\left(
\mathrm{CH}^{(2)}_{s}\left(
\mathfrak{P}^{(2)}
\right)\right)\right)
=
\mathrm{CH}^{(2)}_{s}\left(
\mathfrak{P}^{(2)}
\right).
$$

This completes the proof.
\end{proof}

\begin{restatable}{corollary}{CDThetaCongDCH}
\label{CDThetaCongDCH} Let $s$ be a sort in $S$ and $\mathfrak{P}'^{(2)}$, $\mathfrak{P}^{(2)}$ two second-order paths in $\mathrm{Pth}_{\boldsymbol{\mathcal{A}}^{(2)},s}$. If $(\mathrm{CH}^{(2)}_{s}(\mathfrak{P}'^{(2)}),\mathrm{CH}^{(2)}_{s}(\mathfrak{P}^{(2)}))\in\Theta^{[2]}_{s}$, then $\mathrm{CH}^{(2)}_{s}(\mathfrak{P}'^{(2)})=\mathrm{CH}^{(2)}_{s}(\mathfrak{P}^{(2)})$.
\end{restatable}
\begin{proof}
The following chain of equalities holds
\allowdisplaybreaks
\begin{align*}
\mathrm{CH}^{(2)}_{s}\left(\mathfrak{P}'^{(2)}\right)
&=
\mathrm{CH}^{(2)}_{s}\left(
\mathrm{ip}^{(2,X)@}_{s}\left(
\mathrm{CH}^{(2)}_{s}\left(
\mathfrak{P}'^{(2)}
\right)\right)\right)
\tag{1}
\\&=
\mathrm{CH}^{(2)}_{s}\left(
\mathfrak{P}^{(2)}
\right).\tag{2}
\end{align*}

The first equality follows from Proposition~\ref{PDIpDCH}, note that $\mathrm{ip}^{(2,X)@}_{s}(\mathrm{CH}^{(2)}_{s}(\mathfrak{P}'^{(2)}))$ is a second-order path in $[\mathfrak{P}'^{(2)}]_{s}$; the second equality follows from Corollary~\ref{CDThetaCong}.
\end{proof}

\begin{remark}\label{RDExplFalla} Taking into account Lemmas~\ref{LDWCong} and~\ref{LDThetaCong} and Corollary~\ref{CDThetaCongDCH} we infer that, for every sort $s\in S$, if a term $P\in\mathrm{T}_{\Sigma^{\boldsymbol{\mathcal{A}}^{(2)}}}(X)_{s}$ belongs to $[\mathrm{CH}^{(2)}(\mathfrak{P}^{(2)})]_{\Theta^{[2]}_{s}}$, for some second-order path $\mathfrak{P}^{(2)}\in\mathrm{Pth}_{\boldsymbol{\mathcal{A}}^{(2)},s}$ then $P$ can be understood as an alternative term describing the same second-order path as $\mathfrak{P}^{(2)}$ up to equivalence with respect to $\mathrm{Ker}(\mathrm{CH}^{(2)})$. Therefore, we can think of the class $[\mathrm{CH}^{(2)}_{s}(\mathfrak{P}^{(2)})]_{\Theta^{[2]}_{s}}$ as a class containing all possible terms describing all second-order paths in $[\mathfrak{P}^{(2)}]^{}_{s}$. Given one such term, $P\in[\mathrm{CH}^{(2)}_{s}(\mathfrak{P}^{(2)})]_{\Theta^{[2]}_{s}}$, then $\mathrm{ip}^{(2,X)@}_{s}(P)$ is a second-order path in $[\mathfrak{P}^{(2)}]^{}_{s}$. 
\end{remark}

\chapter{
\texorpdfstring
{The congruence $\Psi^{[1]}$ on $\mathbf{T}_{\Sigma^{\boldsymbol{\mathcal{A}}^{(2)}}}(X)$ and its relatives}
{The congruence Psi}
}\label{S2K}

In this chapter we begin by defining a relation $\Psi^{(1)}$ on $\mathrm{T}_{\Sigma^{\boldsymbol{\mathcal{A}}^{(2)}}}(X)$ that match those terms that represent second-order paths that are $\Upsilon^{(1)}$-related. In fact, the first result of this chapter states that if a pair of terms in $\mathrm{T}_{\Sigma^{\boldsymbol{\mathcal{A}}^{(2)}}}(X)_{s}$ are $\Psi^{(1)}_{s}$-related and, after applying the $\mathrm{ip}^{(2,X)@}$ mapping, they become second-order paths, then these are $\Upsilon^{(1)}_{s}$-related.  We also consider the smallest $\Sigma^{\boldsymbol{\mathcal{A}}^{(2)}}$-congruence containing $\Psi^{(1)}$, denoted by $\Psi^{[1]}$. We will denote by $[\mathbf{T}_{\Sigma^{\boldsymbol{\mathcal{A}}^{(2)}}}(X)]_{\Psi^{[1]}}$ the partial $\Sigma^{\boldsymbol{\mathcal{A}}^{(2)}}$-algebra on the quotient. The $\Sigma^{\boldsymbol{\mathcal{A}}^{(2)}}$-congruence $\Psi^{[1]}$ also satisfies that, for a pair of terms in $\mathrm{T}_{\Sigma^{\boldsymbol{\mathcal{A}}^{(2)}}}(X)_{s}$ if they are $\Psi^{[1]}_{s}$-related and, after applying the $\mathrm{ip}^{(2,X)@}$ mapping, they become second-order paths, then these are $\Upsilon^{[1]}_{s}$-related.


The content of this chapter is unlike anything that has been worked on in the first part. This chapter attempts to reproduce the results of Chapter~\ref{S2F}, where we introduced the $\Sigma^{\boldsymbol{\mathcal{A}}^{(2)}}$-congruence $\Upsilon^{[1]}$ on second-order paths.

We first introduce the relation $\Psi^{(1)}$. In it two terms will be related if they correspond to the terms resulting from the application of the second-order Curry-Howard mapping on paths that are $\Upsilon^{(1)}$-related in $\mathbf{Pth}_{\boldsymbol{\mathcal{A}}^{(2)}}$. 

We will now define the relation we will be working with in this chapter.

\begin{restatable}{definition}{DDPsi}
\label{DDPsi}
\index{Psi!second-order!$\Psi^{(1)}$}
 We define $\Psi^{(1)}=(\Psi^{(1)}_{s})_{s\in S}$ to be the relation on $\mathrm{T}_{\Sigma^{\boldsymbol{\mathcal{A}}^{(2)}}}(X)$ consisting exactly on the following pair of terms
\begin{enumerate}
\item[(i)] For every sort $s\in S$ and every pair of second-order paths $\mathfrak{P}^{(2)}$ and $\mathfrak{Q}^{(2)}$ in $\mathrm{Pth}_{\boldsymbol{\mathcal{A}}^{(2)},s}$ satisfying that $(\mathfrak{P}^{(2)},\mathfrak{Q}^{(2)})\in\Upsilon^{(1)}_{s}$, then 
\[
\left(
\mathrm{CH}^{(2)}_{s}\left(\mathfrak{P}^{(2)}\right),
\mathrm{CH}^{(2)}_{s}\left(\mathfrak{Q}^{(2)}\right)
\right)
\in\Psi^{(1)}_{s}.
\]
\end{enumerate} 
This completes the definition of $\Psi^{(1)}$.
\end{restatable}

Next proposition states that the pairs defining the relation $\Psi^{(1)}$ introduced in Definition~\ref{DDPsi}, when interpreted as a second-order paths, by means of $\mathrm{ip}^{(2,X)@}$ mapping retrieve second-order paths that are related for the relation $\Upsilon^{(1)}$ on $\mathbf{Pth}_{\boldsymbol{\mathcal{A}}^{(2)}}$ introduced in Definition~\ref{DDUps}.

\begin{restatable}{proposition}{PDPsiIp}
\label{PDPsiIp} Let $s$ be a sort in $S$ and $P,Q$ terms in $\mathrm{T}_{\Sigma^{\boldsymbol{\mathcal{A}}^{(2)}}}(X)_{s}$ such that $(P,Q)\in\Psi^{(1)}_{s}$, then
\begin{itemize}
\item[(i)] $\mathrm{ip}^{(2,X)@}_{s}(P)\in\mathrm{Pth}_{\boldsymbol{\mathcal{A}}^{(2)},s}$ if, and only if, $\mathrm{ip}^{(2,X)@}_{s}(Q)\in\mathrm{Pth}_{\boldsymbol{\mathcal{A}}^{(2)},s}$;
\item[(ii)] If $\mathrm{ip}^{(2,X)@}_{s}(P)\in\mathrm{Pth}_{\boldsymbol{\mathcal{A}}^{(2)},s}$ or  $\mathrm{ip}^{(2,X)@}_{s}(Q)\in\mathrm{Pth}_{\boldsymbol{\mathcal{A}}^{(2)},s}$ is a second-order path in $\mathrm{Pth}_{\boldsymbol{\mathcal{A}}^{(2)},s}$ then
$
(
\mathrm{ip}^{(2,X)@}_{s}(
P), 
\mathrm{ip}^{(2,X)@}_{s}(
Q))\in \Upsilon^{(1)}_{s}.
$
\end{itemize}
\end{restatable}
\begin{proof}
Let us recall the definition of the relation $\Psi^{(1)}_{s}$ introduced in Definition~\ref{DDPsi}. For every sort $s\in S$ and every pair of second-order paths $\mathfrak{P}^{(2)}$ and $\mathfrak{Q}^{(2)}$ in $\mathrm{Pth}_{\boldsymbol{\mathcal{A}}^{(2)},s}$ satisfying that $(\mathfrak{P}^{(2)},\mathfrak{Q}^{(2)})\in\Upsilon^{(1)}_{s}$, then 
\[
\left(
\mathrm{CH}^{(2)}_{s}\left(\mathfrak{P}^{(2)}\right),
\mathrm{CH}^{(2)}_{s}\left(\mathfrak{Q}^{(2)}\right)
\right)
\in\Psi^{(1)}_{s}.
\]

According to Proposition~\ref{PDIpDCH} $\mathrm{ip}^{(2,X)@}_{s}(
\mathrm{CH}^{(2)}_{s}(
\mathfrak{P}^{(2)}
))$ and $\mathrm{ip}^{(2,X)@}_{s}(
\mathrm{CH}^{(2)}_{s}(
\mathfrak{Q}^{(2)}
))$ are second-order paths. 

According to Proposition~\ref{PDIpUps} the following pairs of second-order paths are $\Upsilon^{(1)}_{s}$-related
\[
\left(
\mathrm{ip}^{(2,X)@}_{s}\left(
\mathrm{CH}^{(2)}_{s}\left(
\mathfrak{P}^{(2)}
\right)\right),
\mathrm{ip}^{(2,X)@}_{s}\left(
\mathrm{CH}^{(2)}_{s}\left(
\mathfrak{Q}^{(2)}
\right)\right)
\right)
\in\Upsilon^{(1)}_{s}.
\]

This completes the proof.
\end{proof}

\section{
\texorpdfstring
{The congruence $\Psi^{[1]}$ on  $\mathbf{T}_{\Sigma^{\boldsymbol{\mathcal{A}}^{(2)}}}(X)$}
{The congruence Psi}
}

We next show that the properties that hold for the relation $\Psi^{(1)}$, set up in Definition~\ref{DDPsi} also hold for the smallest $\Sigma^{\boldsymbol{\mathcal{A}}^{(2)}}$-congruence generated by it. In this regard, we recover the notation introduced in Definition~\ref{DCongOpInt}.

\begin{restatable}{definition}{DDPsiCong}
\label{DDPsiCong} 
\index{Psi!second-order!$\Psi^{[1]}$}
We denote by $\Psi^{[1]}$ the smallest $\Sigma^{\boldsymbol{\mathcal{A}}^{(2)}}$-congruence on  $\mathbf{T}_{\Sigma^{\boldsymbol{\mathcal{A}}^{(2)}}}(X)$ containing $\Psi^{(1)}$ i.e.,
$$\Psi^{[1]}
=\mathrm{Cg}_{\mathbf{T}_{\Sigma^{\boldsymbol{\mathcal{A}}^{(2)}}}(X)}
\left(\Psi^{(1)}
\right).
$$
\end{restatable}

\begin{remark}\label{RDPsiCong} Let $\mathrm{C}_{\Sigma^{\boldsymbol{\mathcal{A}}^{(2)}}}$ stand for the operator $\mathrm{C}_{\mathbf{T}_{\Sigma^{\boldsymbol{\mathcal{A}}^{(2)}}}(X)}$ on $\mathrm{T}_{\Sigma^{\boldsymbol{\mathcal{A}}^{(2)}}}(X)\times \mathrm{T}_{\Sigma^{\boldsymbol{\mathcal{A}}^{(2)}}}(X)$ (see Definition~\ref{DCongOpC} for the general case). Let us also recall, from Definition~\ref{DCongOpC}, that if $\Phi\subseteq \mathrm{T}_{\Sigma^{\boldsymbol{\mathcal{A}}^{(2)}}}(X)\times \mathrm{T}_{\Sigma^{\boldsymbol{\mathcal{A}}^{(2)}}}(X)$, then 
$$
\mathrm{C}_{\Sigma^{\boldsymbol{\mathcal{A}}^{(2)}}}(\Phi)
=
(\Phi\circ\Phi)
\cup
\left(
\bigcup_{\gamma\in\Sigma^{\boldsymbol{\mathcal{A}}^{(2)}}_{\neq\lambda,s}}
\gamma^{\mathbf{T}_{\Sigma^{\boldsymbol{\mathcal{A}}^{(2)}}}(X)}
\times
\gamma^{\mathbf{T}_{\Sigma^{\boldsymbol{\mathcal{A}}^{(2)}}}(X)}
\left[
\Phi_{\mathrm{ar}(\gamma)}
\right]
\right)_{s\in S}.
$$
Moreover, for the family $(\mathrm{C}^{n}_{\Sigma^{\boldsymbol{\mathcal{A}}^{(2)}}}(\Psi^{(1)}))_{n\in\mathbb{N}}$ in $\mathrm{Sub}(\mathrm{T}_{\Sigma^{\boldsymbol{\mathcal{A}}^{(2)}}}(X)^{2}
)
$, defined recursively as follows:
\allowdisplaybreaks
\begin{align*}
\mathrm{C}^{0}_{\Sigma^{\boldsymbol{\mathcal{A}}^{(2)}}}
\left(
\Psi^{(1)}
\right)
&=\left(\Psi^{(1)}\right)
\cup
\left(
\Psi^{(1)}
\right)^{-1}
\cup\Delta_{\mathrm{T}_{\Sigma^{\boldsymbol{\mathcal{A}}^{(2)}}}(X)},
\\
\mathrm{C}^{n+1}_{\Sigma^{\boldsymbol{\mathcal{A}}^{(2)}}}\left(
\Psi^{(1)}
\right)
&=\mathrm{C}_{\Sigma^{\boldsymbol{\mathcal{A}}^{(2)}}}\left(
\mathrm{C}^{n+1}_{\Sigma^{\boldsymbol{\mathcal{A}}^{(2)}}}\left(
\Psi^{(1)}
\right)\right),\,n\geq 0,
\end{align*}
we have, by Proposition~\ref{PCongOpC}, that 
$$\mathrm{C}^{\omega}_{\Sigma^{\boldsymbol{\mathcal{A}}^{(2)}}}\left(
\Psi^{(1)}
\right)
=
\bigcup_{n\in\mathbb{N}}\mathrm{C}^{n}_{\Sigma^{\boldsymbol{\mathcal{A}}^{(2)}}}\left(
\Psi^{(1)}\right)=\Psi^{[1]},$$ is equal to the smallest $\Sigma^{\boldsymbol{\mathcal{A}}^{(2)}}$-congruence on $\mathbf{T}_{\Sigma^{\boldsymbol{\mathcal{A}}^{(2)}}}(X)$ containing $\Psi^{(1)}$.
\end{remark}

\begin{restatable}{proposition}{PDPsiCongDCatAlg}
\label{PDPsiCongDCatAlg}
\index{terms!second-order!$\mathrm{T}_{\Sigma^{\boldsymbol{\mathcal{A}}^{(2)}}}(X)/{\Psi^{[1]}}$}
\index{terms!second-order!$[P]_{\Psi^{[1]}_{s}}$}
 The $S$-sorted set $\mathrm{T}_{\Sigma^{\boldsymbol{\mathcal{A}}^{(2)}}}(X)/{\Psi^{[1]}}$ is equipped, in a natural way, with a structure of many-sorted $\Sigma^{\boldsymbol{\mathcal{A}}^{(2)}}$-algebra. We denote by $\mathbf{T}_{\Sigma^{\boldsymbol{\mathcal{A}}^{(2)}}}(X)/{\Psi^{[1]}}$ the corresponding $\Sigma^{\boldsymbol{\mathcal{A}}^{(2)}}$-algebra.

\index{projection!second-order!$\mathrm{pr}^{\Psi^{[1]}}$}
We define the projection $\mathrm{pr}^{\Psi^{[1]}}$ from $\mathrm{T}_{\Sigma^{\boldsymbol{\mathcal{A}}^{(2)}}}(X)$ to $\mathrm{T}_{\Sigma^{\boldsymbol{\mathcal{A}}^{(2)}}}(X)/{\Psi^{[1]}}$ to be the many-sorted mapping that, for every sort $s\in S$, maps a term $P$ in $\mathrm{T}_{\Sigma^{\boldsymbol{\mathcal{A}}^{(2)}}}(X)$ to $[P]_{\Psi^{[1]}_{s}}$, its equivalence class under the $\Sigma^{\boldsymbol{\mathcal{A}}^{(2)}}$-congruence $\Psi^{[1]}$.

Note that $\mathrm{pr}^{\Psi^{[1]}}$, that is,
$$
\mathrm{pr}^{\Psi^{[1]}}\colon
\mathbf{T}_{\Sigma^{\boldsymbol{\mathcal{A}}^{(2)}}}(X)
\mor
\mathbf{T}_{\Sigma^{\boldsymbol{\mathcal{A}}^{(2)}}}(X)/{\Psi^{[1]}}
$$
is a surjective $\Sigma^{\boldsymbol{\mathcal{A}}^{(2)}}$-homomorphism from $\mathbf{T}_{\Sigma^{\boldsymbol{\mathcal{A}}^{(2)}}}(X)$ to $\mathbf{T}_{\Sigma^{\boldsymbol{\mathcal{A}}^{(2)}}}(X)/{\Psi^{[1]}}$.
\end{restatable}

We next introduce the reductions of the $\Sigma^{\boldsymbol{\mathcal{A}}^{(2)}}$-algebra $\mathbf{T}_{\Sigma^{\boldsymbol{\mathcal{A}}^{(2)}}}(X)/{\Psi^{[1]}}$ to layers $0$ and $1$.

\begin{definition} We will denote by $\mathbf{T}^{(1,2)}_{\Sigma^{\boldsymbol{\mathcal{A}}^{(2)}}}(X)/{\Psi^{[1]}}$ the $\Sigma^{\boldsymbol{\mathcal{A}}}$-algebra 
$$\left(\mathbf{in}^{\Sigma,(1,2)}
\right)\left(\mathbf{T}_{\Sigma^{\boldsymbol{\mathcal{A}}^{(2)}}}(X)/{\Psi^{[1]}}\right).$$

We will call $\mathbf{T}^{(1,2)}_{\Sigma^{\boldsymbol{\mathcal{A}}^{(2)}}}(X)/{\Psi^{[1]}}$ the \emph{$\Sigma^{\boldsymbol{\mathcal{A}}}$-reduct} of the many-sorted $\Sigma^{\boldsymbol{\mathcal{A}}^{(2)}}$-algebra $\mathbf{T}_{\Sigma^{\boldsymbol{\mathcal{A}}^{(2)}}}(X)/{\Psi^{[1]}}$. In this regard see Definition~\ref{DDURed} and Remark~\ref{RDURed}.

We will denote by $\mathbf{T}^{(0,2)}_{\Sigma^{\boldsymbol{\mathcal{A}}^{(2)}}}(X)/{\Psi^{[1]}}$ the $\Sigma$-algebra 
$$\left(\mathbf{in}^{\Sigma,(0,2)}
\right)\left(\mathbf{T}_{\Sigma^{\boldsymbol{\mathcal{A}}^{(2)}}}(X)/{\Psi^{[1]}}\right).$$

We will call $\mathbf{T}^{(0,2)}_{\Sigma^{\boldsymbol{\mathcal{A}}^{(2)}}}(X)/{\Psi^{[1]}}$ the \emph{$\Sigma$-reduct} of the many-sorted $\Sigma^{\boldsymbol{\mathcal{A}}^{(2)}}$-algebra $\mathbf{T}_{\Sigma^{\boldsymbol{\mathcal{A}}^{(2)}}}(X)/{\Psi^{[1]}}$. In this regard see Definition~\ref{DDZRed} and Remark~\ref{RDZRed}.
\end{definition}

The following lemma states that, for every sort $s\in S$, if two terms in $\mathrm{T}_{\Sigma^{\boldsymbol{\mathcal{A}}^{(2)}}}(X)_{s}$ are $\Psi^{[1]}_{s}$-related and one of them, when mapped under $\mathrm{ip}^{(2,X)@}_{s}$, retrieves a second-order path, then the other term has a similar behaviour. Moreover, if the just described situation happens, then these two second-order paths will be $\Upsilon^{[1]}_{s}$-related in $\mathrm{Pth}_{\boldsymbol{\mathcal{A}}^{(2)},s}$. The interested reader can see how the proof of the following result resembles the proof of the Lemma~\ref{LDThetaCong}.

\begin{restatable}{lemma}{LDPsiCong}
\label{LDPsiCong} Let $s$ be a sort in $S$ and $P,Q$ terms in $\mathrm{T}_{\Sigma^{\boldsymbol{\mathcal{A}}^{(2)}}}(X)_{s}$ such that $(P,Q)\in\Psi^{[1]}_{s}$, then
\begin{itemize}
\item[(i)] $\mathrm{ip}^{(2,X)@}_{s}(P)\in\mathrm{Pth}_{\boldsymbol{\mathcal{A}}^{(2)},s}$ if, and only if, $\mathrm{ip}^{(2,X)@}_{s}(Q)\in\mathrm{Pth}_{\boldsymbol{\mathcal{A}}^{(2)},s}$;
\item[(ii)] If $\mathrm{ip}^{(2,X)@}_{s}(P)\in\mathrm{Pth}_{\boldsymbol{\mathcal{A}}^{(2)},s}$ or  $\mathrm{ip}^{(2,X)@}_{s}(Q)\in\mathrm{Pth}_{\boldsymbol{\mathcal{A}}^{(2)},s}$ is a second-order path in $\mathrm{Pth}_{\boldsymbol{\mathcal{A}}^{(2)},s}$ then
$
(
\mathrm{ip}^{(2,X)@}_{s}(
P), 
\mathrm{ip}^{(2,X)@}_{s}(
Q))\in \Upsilon^{[1]}_{s}.
$
\end{itemize}
\end{restatable}
\begin{proof}
We recall from Remark~\ref{RDPsiCong} that 
$$
\Psi^{[1]}=
\bigcup_{n\in\mathbb{N}}\mathrm{C}^{n}_{\Sigma^{\boldsymbol{\mathcal{A}}^{(2)}}}
\left(
\Psi^{(1)}\right).
$$

We prove the statement by induction on $n\in\mathbb{N}$.

\textsf{Base step of the induction.}

Let us recall from Remark~\ref{RDPsiCong} that 
\[
\mathrm{C}^{0}_{\Sigma^{\boldsymbol{\mathcal{A}}^{(2)}}}\left(\Psi^{(1)}\right)=
\left(\Psi^{(1)}\right)
\cup
\left(
\Psi^{(1)}
\right)^{-1}
\cup
\Delta_{\mathrm{T}_{\Sigma^{\boldsymbol{\mathcal{A}}^{(2)}}}(X)}.
\]

The statement trivially holds for a pair $(P,Q)$ in $\Delta_{\mathrm{T}_{\Sigma^{\boldsymbol{\mathcal{A}}^{(2)}}}(X)_{s}}$.

If the pair $(P,Q)$ is in $\Psi^{(1)}_{s}$, following Proposition~\ref{PDPsiIp}, we have that both $\mathrm{ip}^{(2,X)@}_{s}(P)$ and $\mathrm{ip}^{(2,X)@}_{s}(Q)$ second-order paths that are $\Upsilon^{(1)}_{s}$-related. The same applies if $(P,Q)$ is a pair in  $(\Psi^{(1)}_{s})^{-1}$. 

This completes the base case.

\textsf{Inductive step of the induction.}

Assume the statement holds for $n\in\mathbb{N}$, i.e., for every sort $s\in S$ and every pair of terms $(P,Q)$ in $\mathrm{T}_{\Sigma^{\boldsymbol{\mathcal{A}}^{(2)}}}(X)_{s}$ such that $(P,Q)$ in $\mathrm{C}^{n}_{\Sigma^{\boldsymbol{\mathcal{A}}^{(2)}}}(\Psi^{(1)})_{s}$ then
\begin{itemize}
\item[(i)] $\mathrm{ip}^{(2,X)@}_{s}(P)\in\mathrm{Pth}_{\boldsymbol{\mathcal{A}}^{(2)},s}$ if, and only if, $\mathrm{ip}^{(2,X)@}_{s}(Q)\in\mathrm{Pth}_{\boldsymbol{\mathcal{A}}^{(2)},s}$;
\item[(ii)] If $\mathrm{ip}^{(2,X)@}_{s}(P)\in\mathrm{Pth}_{\boldsymbol{\mathcal{A}}^{(2)},s}$ or  $\mathrm{ip}^{(2,X)@}_{s}(Q)\in\mathrm{Pth}_{\boldsymbol{\mathcal{A}}^{(2)},s}$ is a second-order path in $\mathrm{Pth}_{\boldsymbol{\mathcal{A}}^{(2)},s}$ then
$
(
\mathrm{ip}^{(2,X)@}_{s}(
P), 
\mathrm{ip}^{(2,X)@}_{s}(
Q))\in \Upsilon^{[1]}_{s}.
$
\end{itemize}

We prove the statement for $n+1$. Let $s$ be a sort in $S$ and let $(P,Q)$ be a pair of terms in $\mathrm{T}_{\Sigma^{\boldsymbol{\mathcal{A}}^{(2)}}}(X)_{s}$ such that $(P,Q)$ is a pair in $\mathrm{C}_{\Sigma^{\boldsymbol{\mathcal{A}}^{(2)}}}^{n+1}(\Psi^{(1)})_{s}$. Let us recall from Remark~\ref{RDPsiCong} that 
\begin{multline*}
\mathrm{C}^{n+1}_{\Sigma^{\boldsymbol{\mathcal{A}}^{(2)}}}\left(
\Psi^{(1)}
\right)_{s}
=
\left(
\mathrm{C}^{n}_{\Sigma^{\boldsymbol{\mathcal{A}}^{(2)}}}\left(
\Psi^{(1)}
\right)_{s}
\circ
\mathrm{C}^{n}_{\Sigma^{\boldsymbol{\mathcal{A}}^{(2)}}}\left(
\Psi^{(1)}
\right)_{s}
\right)
\cup
\\
\left(
\bigcup_{\gamma\in\Sigma^{\boldsymbol{\mathcal{A}}^{(2)}}_{\neq\lambda,s}}
\gamma^{\mathbf{T}_{\Sigma^{\boldsymbol{\mathcal{A}}^{(2)}}}(X)}
\times
\gamma^{\mathbf{T}_{\Sigma^{\boldsymbol{\mathcal{A}}^{(2)}}}(X)}
\left[
\mathrm{C}^{n}_{\Sigma^{\boldsymbol{\mathcal{A}}^{(2)}}}\left(
\Psi^{(1)}
\right)_{\mathrm{ar}(\gamma)}
\right]
\right).
\end{multline*}

Then either (1), $(P,Q)$ is in $\mathrm{C}^{n}_{\Sigma^{\boldsymbol{\mathcal{A}}^{(2)}}}(
\Psi^{(1)}
)_{s}
\circ
\mathrm{C}^{n}_{\Sigma^{\boldsymbol{\mathcal{A}}^{(2)}}}(
\Psi^{(1)}
)_{s}$ or (2), $(P,Q)$ is in $\gamma^{\mathbf{T}_{\Sigma^{\boldsymbol{\mathcal{A}}^{(2)}}}(X)}
\times
\gamma^{\mathbf{T}_{\Sigma^{\boldsymbol{\mathcal{A}}^{(2)}}}(X)}
[
\mathrm{C}^{n}_{\Sigma^{\boldsymbol{\mathcal{A}}^{(2)}}}(
\Psi^{(1)}
)_{\mathrm{ar}(\gamma)}
]$ for some operation symbol $\gamma\in\Sigma^{\boldsymbol{\mathcal{A}}^{(2)}}_{\neq\lambda,s}$.

If (1), then there exists a term $R\in\mathrm{T}_{\Sigma^{\boldsymbol{\mathcal{A}}^{(2)}}}(X)_{s}$ for which $(P,R)$ and $(R,Q)$ belong to $\mathrm{C}^{n}_{\Sigma^{\boldsymbol{\mathcal{A}}^{(2)}}}(\Psi^{(1)})_{s}$. Then by induction we have that 
\[
\mathrm{ip}^{(2,X)@}_{s}(P)\in\mathrm{Pth}_{\boldsymbol{\mathcal{A}}^{(2)},s}
\Leftrightarrow
\mathrm{ip}^{(2,X)@}_{s}(R)\in\mathrm{Pth}_{\boldsymbol{\mathcal{A}}^{(2)},s}
\Leftrightarrow
\mathrm{ip}^{(2,X)@}_{s}(Q)\in\mathrm{Pth}_{\boldsymbol{\mathcal{A}}^{(2)},s}.
\]

Moreover, in case on of the previous elements is a second-order path in $\mathrm{Pth}_{\boldsymbol{\mathcal{A}}^{(2)}}$, we  have, by induction, that 
\allowdisplaybreaks
\begin{align*}
\left(
\mathrm{ip}^{(2,X)@}_{s}\left(P\right),
\mathrm{ip}^{(2,X)@}_{s}\left(R\right)
\right)
&\in\Upsilon^{[1]}_{s};
&&
\left(
\mathrm{ip}^{(2,X)@}_{s}\left(R\right),
\mathrm{ip}^{(2,X)@}_{s}\left(Q\right)
\right)
&\in\Upsilon^{[1]}_{s}.
\end{align*}
Since $\Upsilon^{[1]}$ is an equivalence relation on $\mathbf{Pth}_{\boldsymbol{\mathcal{A}}^{(2)}}$, we conclude that 
\[
\left(
\mathrm{ip}^{(2,X)@}_{s}\left(P\right),
\mathrm{ip}^{(2,X)@}_{s}\left(Q\right)
\right)
\in\Upsilon^{[1]}_{s}.
\]

If (2), then there exists a unique word $\mathbf{s}\in S^{\star}-\{\lambda\}$, a unique operation symbol $\gamma\in \Sigma^{\boldsymbol{\mathcal{A}}^{(2)}}_{\mathbf{s},s}$ and a unique family of pairs $((P_{j},Q_{j}))_{j\in\bb{\mathbf{s}}}$ in $\mathrm{C}^{n}_{\Sigma^{\boldsymbol{\mathcal{A}}^{(2)}}}(\Psi^{(1)})_{\mathbf{s}}$ for which
\[
(P,Q)=\left(
\gamma^{\mathbf{T}_{\Sigma^{\boldsymbol{\mathcal{A}}^{(2)}}}(X)}
\left(\left(
P_{j}\right)_{j\in\bb{\mathbf{s}}}\right),
\gamma^{\mathbf{T}_{\Sigma^{\boldsymbol{\mathcal{A}}^{(2)}}}(X)}
\left(\left(
Q_{j}\right)_{j\in\bb{\mathbf{s}}}\right)
\right).
\]

We will distinguish the following cases according to the different possibilities for the operation symbol $\gamma\in\Sigma^{\boldsymbol{\mathcal{A}}^{(2)}}_{\mathbf{s},s}$. Note that either $\gamma$ is an operation symbol $\sigma\in\Sigma_{\mathbf{s},s}$, the operation symbol of $0$-source, $\mathrm{sc}_{s}^{0}$,  the operation symbol of $0$-target, $\mathrm{tg}_{s}^{0}$,  the operation symbol of $0$-composition, $\circ^{0}_{s}$,  the operation symbol of $1$-source, $\mathrm{sc}_{s}^{1}$,  the operation symbol of $1$-target, $\mathrm{tg}_{s}^{1}$, or the operation symbol of $1$-composition, $\circ^{1}_{s}$.

\textsf{$\gamma$ is an operation symbol $\sigma\in\Sigma_{\mathbf{s},s}$.}

Let $\mathbf{s}$ be a word in $S^{\star}-\{\lambda\}$, let $\sigma$ be an operation symbol in $\Sigma_{\mathbf{s},s}$ and let $((P_{j},Q_{j}))_{j\in\bb{\mathbf{s}}}$ be a family of pairs in $\mathrm{C}^{n}_{\Sigma^{\boldsymbol{\mathcal{A}}^{(2)}}}(\Psi^{(1)})_{\mathbf{s}}$ for which 
\[
(P,Q)=
\left(
\sigma^{\mathbf{T}_{\Sigma^{\boldsymbol{\mathcal{A}}^{(2)}}}(X)}
\left(\left(
P_{j}\right)_{j\in\bb{\mathbf{s}}}\right),
\sigma^{\mathbf{T}_{\Sigma^{\boldsymbol{\mathcal{A}}^{(2)}}}(X)}
\left(\left(
Q_{j}\right)_{j\in\bb{\mathbf{s}}}\right)
\right).
\]

Note that the following chain of equivalences holds
\allowdisplaybreaks
\begin{align*}
\mathrm{ip}^{(2,X)@}_{s}(P)\in\mathrm{Pth}_{\boldsymbol{\mathcal{A}}^{(2)},s}
&\Leftrightarrow
\mbox{for every }j\in\bb{\mathbf{s}},\,\mathrm{ip}^{(2,X)@}_{s_{j}}(P_{j})\in\mathrm{Pth}_{\boldsymbol{\mathcal{A}}^{(2)},s_{j}}
\tag{1}
\\&\Leftrightarrow
\mbox{for every }j\in\bb{\mathbf{s}},\,\mathrm{ip}^{(2,X)@}_{s_{j}}(Q_{j})\in\mathrm{Pth}_{\boldsymbol{\mathcal{A}}^{(2)},s_{j}}
\tag{2}
\\&\Leftrightarrow
\mathrm{ip}^{(2,X)@}_{s}(Q)\in\mathrm{Pth}_{\boldsymbol{\mathcal{A}}^{(2)},s}.
\tag{3}
\end{align*}

In the just stated chain of equivalences, the first equivalence follows from left to right by Lemma~\ref{LDThetaCongSub} and from right to left because $\sigma$ is a total operation in $\mathbf{Pth}_{\boldsymbol{\mathcal{A}}^{(2)}}$, according to Proposition~\ref{PDPthAlg}, and $\mathrm{ip}^{(2,X)@}$ is a $\Sigma^{\boldsymbol{\mathcal{A}}^{(2)}}$-homomorphism, according to Definition~\ref{DDIp}; the second equivalence follows by induction, finally the third equivalence follows from left to right because $\sigma$ is a total operation in $\mathbf{Pth}_{\boldsymbol{\mathcal{A}}^{(2)}}$, according to Proposition~\ref{PDPthAlg},  and $\mathrm{ip}^{(2,X)@}$ is a $\Sigma^{\boldsymbol{\mathcal{A}}^{(2)}}$-homomorphism, according to Definition~\ref{DDIp}, and from right to left by Lemma~\ref{LDThetaCongSub}.

Assume, without loss of generality, that $\mathrm{ip}^{(2,X)@}_{s}(P)$ is a second-order path in $\mathrm{Pth}_{\boldsymbol{\mathcal{A}}^{(2)},s}$. As we have seen before, this is the case exactly when, for every $j\in\bb{\mathbf{s}}$, $\mathrm{ip}^{(2,X)@}_{s_{j}}(P_{j})$ is a second-order path in $\mathrm{Pth}_{\boldsymbol{\mathcal{A}}^{(2)},s_{j}}$. By induction, we also have that, for every $j\in\bb{\mathbf{s}}$, $\mathrm{ip}^{(2,X)@}_{s_{j}}(Q_{j})$ is a second-order path in $\mathrm{Pth}_{\boldsymbol{\mathcal{A}}^{(2)},s_{j}}$. Moreover, by induction, the following statement holds
\[
\left(
\mathrm{ip}^{(2,X)@}_{s_{j}}\left(
P_{j}
\right),
\mathrm{ip}^{(2,X)@}_{s_{j}}\left(
Q_{j}
\right)\right)\in\Upsilon^{[1]}_{s_{j}}.
\]

Note that the following chain of equalities holds
\allowdisplaybreaks
\begin{align*}
\mathrm{pr}^{\Upsilon^{[1]}}_{s}\left(
\mathrm{ip}^{(2,X)@}_{s}\left(
P
\right)
\right)
&=
\mathrm{pr}^{\Upsilon^{[1]}}_{s}\left(
\mathrm{ip}^{(2,X)@}_{s}\left(
\sigma^{\mathbf{T}_{\Sigma^{\boldsymbol{\mathcal{A}}^{(2)}}}(X)}
\left(
\left(
P_{j}
\right)_{j\in\bb{\mathbf{s}}}
\right)
\right)
\right)
\tag{1}
\\&=
\mathrm{pr}^{\Upsilon^{[1]}}_{s}\left(
\sigma^{\mathbf{Pth}_{\boldsymbol{\mathcal{A}}^{(2)}}}
\left(
\left(
\mathrm{ip}^{(2,X)@}_{s_{j}}\left(
P_{j}
\right)
\right)_{j\in\bb{\mathbf{s}}}
\right)
\right)
\tag{2}
\\&=
\mathrm{pr}^{\Upsilon^{[1]}}_{s}\left(
\sigma^{\mathbf{Pth}_{\boldsymbol{\mathcal{A}}^{(2)}}}
\left(
\left(
\mathrm{ip}^{(2,X)@}_{s_{j}}\left(
Q_{j}
\right)
\right)_{j\in\bb{\mathbf{s}}}
\right)
\right)
\tag{3}
\\&=
\mathrm{pr}^{\Upsilon^{[1]}}_{s}\left(
\mathrm{ip}^{(2,X)@}_{s}\left(
\sigma^{\mathbf{T}_{\Sigma^{\boldsymbol{\mathcal{A}}^{(2)}}}(X)}
\left(
\left(
Q_{j}
\right)_{j\in\bb{\mathbf{s}}}
\right)
\right)
\right)
\tag{4}
\\&=
\mathrm{pr}^{\Upsilon^{[1]}}_{s}\left(
\mathrm{ip}^{(2,X)@}_{s}\left(
Q
\right)
\right).
\tag{5}
\end{align*}

In the just stated chain of equalities, the first equality unravels the description of $P$; the second equivalence follows from the fact that $\sigma$ is a total operation in $\mathbf{Pth}_{\boldsymbol{\mathcal{A}}^{(2)}}$, according to Proposition~\ref{PDPthAlg},  and $\mathrm{ip}^{(2,X)@}$ is a $\Sigma^{\boldsymbol{\mathcal{A}}^{(2)}}$-homomorphism, according to Definition~\ref{DDIp}, and the fact that we are assuming that $\mathrm{ip}^{(2,X)@}_{s}(P)$ is a second-order path in $\mathrm{Pth}_{\boldsymbol{\mathcal{A}}^{(2)},s}$; the third equality follows by induction and from the fact that, according to Definition~\ref{DDUpsCong}, $\Upsilon^{[1]}$ is a $\Sigma^{\boldsymbol{\mathcal{A}}^{(2)}}$-congruence; the fourth equality follows from the fact that $\sigma$ is a total operation in $\mathbf{Pth}_{\boldsymbol{\mathcal{A}}^{(2)}}$, according to Proposition~\ref{PDPthAlg},  and $\mathrm{ip}^{(2,X)@}$ is a $\Sigma^{\boldsymbol{\mathcal{A}}^{(2)}}$-homomorphism, according to Definition~\ref{DDIp}, and the fact that we are assuming that, for every $j\in\bb{\mathbf{s}}$, $\mathrm{ip}^{(2,X)@}_{s_{j}}(Q_{j})$ is a second-order path in $\mathrm{Pth}_{\boldsymbol{\mathcal{A}}^{(2)},s_{j}}$; finally, the last equality recovers the description of $Q$.

This completes the case $\sigma\in\Sigma_{\mathbf{s},s}$.

\textsf{$\gamma$ is the $0$-source operation symbol.}

Let $\gamma$ be the operation symbol $\mathrm{sc}^{0}_{s}$ in $\Sigma^{\boldsymbol{\mathcal{A}}}_{s,s}$ and let $(P',Q')$ be a family of pairs in $\mathrm{C}^{n}_{\Sigma^{\boldsymbol{\mathcal{A}}^{(2)}}}(\Psi^{(1)})_{s}$ for which 
$$
(P,Q)=
\left(
\mathrm{sc}^{0\mathbf{T}_{\Sigma^{\boldsymbol{\mathcal{A}}^{(2)}}}(X)}_{s}\left(
P'
\right)
,
\mathrm{sc}^{0\mathbf{T}_{\Sigma^{\boldsymbol{\mathcal{A}}^{(2)}}}(X)}_{s}\left(
Q'
\right)
\right).
$$

Note that the following chain of equivalences holds
\allowdisplaybreaks
\begin{align*}
\mathrm{ip}^{(2,X)@}_{s}(P)\in\mathrm{Pth}_{\boldsymbol{\mathcal{A}}^{(2)},s}
&\Leftrightarrow
\mathrm{ip}^{(2,X)@}_{s}(P')\in\mathrm{Pth}_{\boldsymbol{\mathcal{A}}^{(2)},s}
\tag{1}
\\&\Leftrightarrow
\mathrm{ip}^{(2,X)@}_{s}(Q')\in\mathrm{Pth}_{\boldsymbol{\mathcal{A}}^{(2)},s}
\tag{2}
\\&\Leftrightarrow
\mathrm{ip}^{(2,X)@}_{s}(Q)\in\mathrm{Pth}_{\boldsymbol{\mathcal{A}}^{(2)},s}.
\tag{3}
\end{align*}

In the just stated chain of equivalences, the first equivalence follows from left to right by Lemma~\ref{LDThetaCongSub} and from right to left because $\mathrm{sc}^{0}_{s}$ is a total operation in $\mathbf{Pth}_{\boldsymbol{\mathcal{A}}^{(2)}}$, according to Proposition~\ref{PDPthCatAlg}, and $\mathrm{ip}^{(2,X)@}$ is a $\Sigma^{\boldsymbol{\mathcal{A}}^{(2)}}$-homomorphism, according to Definition~\ref{DDIp}; the second equivalence follows by induction, finally the third equivalence follows from left to right because $\mathrm{sc}^{0}_{s}$ is a total operation in $\mathbf{Pth}_{\boldsymbol{\mathcal{A}}^{(2)}}$, according to Proposition~\ref{PDPthAlg}, and $\mathrm{ip}^{(2,X)@}$ is a $\Sigma^{\boldsymbol{\mathcal{A}}^{(2)}}$-homomorphism, according to Definition~\ref{DDIp}, and from right to left by Lemma~\ref{LDThetaCongSub}.

Assume, without loss of generality, that $\mathrm{ip}^{(2,X)@}_{s}(P)$ is a second-order path in $\mathrm{Pth}_{\boldsymbol{\mathcal{A}}^{(2)},s}$. As we have seen before, this is the case exactly when $\mathrm{ip}^{(2,X)@}_{s}(P')$ is a second-order path in $\mathrm{Pth}_{\boldsymbol{\mathcal{A}}^{(2)},s}$. By induction, we also have that $\mathrm{ip}^{(2,X)@}_{s}(Q')$ is a second-order path in $\mathrm{Pth}_{\boldsymbol{\mathcal{A}}^{(2)},s}$. Moreover, by induction, the following statement holds
\[
\left(
\mathrm{ip}^{(2,X)@}_{s}\left(
P'
\right),
\mathrm{ip}^{(2,X)@}_{s}\left(
Q'
\right)\right)\in\Upsilon^{[1]}_{s}.
\]

Note that the following chain of equalities holds
\allowdisplaybreaks
\begin{align*}
\mathrm{pr}^{\Upsilon^{[1]}}_{s}\left(
\mathrm{ip}^{(2,X)@}_{s}\left(
P
\right)
\right)
&=
\mathrm{pr}^{\Upsilon^{[1]}}_{s}\left(
\mathrm{ip}^{(2,X)@}_{s}\left(
\mathrm{sc}^{0\mathbf{T}_{\Sigma^{\boldsymbol{\mathcal{A}}^{(2)}}}(X)}_{s}
\left(
P'
\right)
\right)
\right)
\tag{1}
\\&=
\mathrm{pr}^{\Upsilon^{[1]}}_{s}\left(
\mathrm{sc}^{0\mathbf{Pth}_{\boldsymbol{\mathcal{A}}^{(2)}}}_{s}
\left(
\mathrm{ip}^{(2,X)@}_{s}\left(
P'
\right)
\right)
\right)
\tag{2}
\\&=
\mathrm{pr}^{\Upsilon^{[1]}}_{s}\left(
\mathrm{sc}^{0\mathbf{Pth}_{\boldsymbol{\mathcal{A}}^{(2)}}}_{s}
\left(
\mathrm{ip}^{(2,X)@}_{s}\left(
Q'
\right)
\right)
\right)
\tag{3}
\\&=
\mathrm{pr}^{\Upsilon^{[1]}}_{s}\left(
\mathrm{ip}^{(2,X)@}_{s}\left(
\mathrm{sc}^{0\mathbf{T}_{\Sigma^{\boldsymbol{\mathcal{A}}^{(2)}}}(X)}_{s}
\left(
Q'
\right)
\right)
\right)
\tag{4}
\\&=
\mathrm{pr}^{\Upsilon^{[1]}}_{s}\left(
\mathrm{ip}^{(2,X)@}_{s}\left(
Q
\right)
\right).
\tag{5}
\end{align*}

In the just stated chain of equalities, the first equality unravels the description of $P$; the second equivalence follows from the fact that $\mathrm{sc}^{0}_{s}$ is a total operation in $\mathbf{Pth}_{\boldsymbol{\mathcal{A}}^{(2)}}$, according to Proposition~\ref{PDPthCatAlg},  and $\mathrm{ip}^{(2,X)@}$ is a $\Sigma^{\boldsymbol{\mathcal{A}}^{(2)}}$-homomorphism, according to Definition~\ref{DDIp}, and the fact that we are assuming that $\mathrm{ip}^{(2,X)@}_{s}(P)$ is a second-order path in $\mathrm{Pth}_{\boldsymbol{\mathcal{A}}^{(2)},s}$; the third equality follows by induction and from the fact that, according to Definition~\ref{DDUpsCong}, $\Upsilon^{[1]}$ is a $\Sigma^{\boldsymbol{\mathcal{A}}^{(2)}}$-congruence; the fourth equality follows from the fact that $\mathrm{sc}^{0}_{s}$ is a total operation in $\mathbf{Pth}_{\boldsymbol{\mathcal{A}}^{(2)}}$, according to Proposition~\ref{PDPthCatAlg},  and $\mathrm{ip}^{(2,X)@}$ is a $\Sigma^{\boldsymbol{\mathcal{A}}^{(2)}}$-homomorphism, according to Definition~\ref{DDIp}, and the fact that we are assuming that $\mathrm{ip}^{(2,X)@}_{s}(Q')$ is a second-order path in $\mathrm{Pth}_{\boldsymbol{\mathcal{A}}^{(2)},s}$; finally, the last equality recovers the description of $Q$.

This completes the case of the $0$-source.

\textsf{$\gamma$ is the $0$-target operation symbol.}

Let $\gamma$ be the operation symbol $\mathrm{tg}^{0}_{s}$ in $\Sigma^{\boldsymbol{\mathcal{A}}}_{s,s}$ and let $(P',Q')$ be a family of pairs in $\mathrm{C}^{n}_{\Sigma^{\boldsymbol{\mathcal{A}}^{(2)}}}(\Psi^{(1)})_{s}$ for which 
$$
(P,Q)=
\left(
\mathrm{tg}^{0\mathbf{T}_{\Sigma^{\boldsymbol{\mathcal{A}}^{(2)}}}(X)}_{s}\left(
P'
\right)
,
\mathrm{tg}^{0\mathbf{T}_{\Sigma^{\boldsymbol{\mathcal{A}}^{(2)}}}(X)}_{s}\left(
Q'
\right)
\right).
$$

Then, the following properties hold
\begin{itemize}
\item[(i)] $\mathrm{ip}^{(2,X)@}_{s}(P)\in\mathrm{Pth}_{\boldsymbol{\mathcal{A}}^{(2)},s}$ if, and only if, $\mathrm{ip}^{(2,X)@}_{s}(Q)\in\mathrm{Pth}_{\boldsymbol{\mathcal{A}}^{(2)},s}$;
\item[(ii)] If $\mathrm{ip}^{(2,X)@}_{s}(P)\in\mathrm{Pth}_{\boldsymbol{\mathcal{A}}^{(2)},s}$ or  $\mathrm{ip}^{(2,X)@}_{s}(Q)\in\mathrm{Pth}_{\boldsymbol{\mathcal{A}}^{(2)},s}$ is a second-order path in $\mathrm{Pth}_{\boldsymbol{\mathcal{A}}^{(2)},s}$ then
$
(
\mathrm{ip}^{(2,X)@}_{s}(
P),
\mathrm{ip}^{(2,X)@}_{s}(
Q))\in\Upsilon^{[1]}_{s}.
$
\end{itemize}

The proof of this case is similar to that of the $0$-source.

This completes the case of the $0$-target.

\textsf{$\gamma$ is the $0$-composition operation symbol.}

Let $\gamma$ be the operation symbol $\circ^{0}_{s}$ in $\Sigma^{\boldsymbol{\mathcal{A}}}_{s,s}$ and let $(P',Q')$ and $(P'',Q'')$ be two families of pairs in $\mathrm{C}^{n}_{\Sigma^{\boldsymbol{\mathcal{A}}^{(2)}}}(\Psi^{(1)})_{s}$ for which 
$$
(P,Q)=
\left(
P''
\circ^{0\mathbf{T}_{\Sigma^{\boldsymbol{\mathcal{A}}^{(2)}}}(X)}_{s}
P'
,
Q''
\circ^{0\mathbf{T}_{\Sigma^{\boldsymbol{\mathcal{A}}^{(2)}}}(X)}_{s}
Q'
\right).
$$

Note that the following chain of equivalences holds
\allowdisplaybreaks
\begin{align*}
\mathrm{ip}^{(2,X)@}_{s}(P)\in\mathrm{Pth}_{\boldsymbol{\mathcal{A}}^{(2)},s}
&\Leftrightarrow
\left\lbrace
\begin{array}{l}
\mathrm{ip}^{(2,X)@}_{s}(P'),\mathrm{ip}^{(2,X)@}_{s}(P'')\in\mathrm{Pth}_{\boldsymbol{\mathcal{A}}^{(2)},s}
\\
\mathrm{sc}^{(0,2)}_{s}(
\mathrm{ip}^{(2,X)@}_{s}(P'')
)
=
\mathrm{tg}^{(0,2)}_{s}(
\mathrm{ip}^{(2,X)@}_{s}(P')
)
\end{array}
\right.
\tag{1}
\\&\Leftrightarrow
\left\lbrace
\begin{array}{l}
\mathrm{ip}^{(2,X)@}_{s}(Q'),\mathrm{ip}^{(2,X)@}_{s}(Q'')\in\mathrm{Pth}_{\boldsymbol{\mathcal{A}}^{(2)},s}
\\
\mathrm{sc}^{(0,2)}_{s}(
\mathrm{ip}^{(2,X)@}_{s}(Q'')
)
=
\mathrm{tg}^{(0,2)}_{s}(
\mathrm{ip}^{(2,X)@}_{s}(Q')
)
\end{array}
\right.
\tag{2}
\\&\Leftrightarrow
\mathrm{ip}^{(2,X)@}_{s}(Q)\in\mathrm{Pth}_{\boldsymbol{\mathcal{A}}^{(2)},s}.
\tag{3}
\end{align*}

In the just stated chain of equivalences, the first equivalence follows from left to right by Lemma~\ref{LDThetaCongSub} and by the description of the $0$-composition operation in $\mathbf{Pth}_{\boldsymbol{\mathcal{A}}^{(2)}}$, according to Proposition~\ref{PDPthCatAlg}, and from right to left because $\circ^{0}_{s}$ is a defined operation in $\mathbf{Pth}_{\boldsymbol{\mathcal{A}}^{(2)}}$, according to Proposition~\ref{PDPthCatAlg}, and $\mathrm{ip}^{(2,X)@}$ is a $\Sigma^{\boldsymbol{\mathcal{A}}^{(2)}}$-homomorphism, according to Definition~\ref{DDIp}; the second equivalence follows by induction. Note also that
\allowdisplaybreaks
\begin{align*}
\left(
\mathrm{ip}^{(2,X)@}_{s}\left(
P'
\right)
,
\mathrm{ip}^{(2,X)@}_{s}\left(
Q'
\right)
\right)&\in\Upsilon^{[1]}_{s};
\\
\left(
\mathrm{ip}^{(2,X)@}_{s}\left(
P''
\right)
,
\mathrm{ip}^{(2,X)@}_{s}\left(
Q''
\right)
\right)&\in\Upsilon^{[1]}_{s}.
\end{align*}
Therefore, according to Corollary~\ref{CDUpsCongDZ}, we have that 
\allowdisplaybreaks
\begin{align*}
\mathrm{sc}^{(0,2)}_{s}\left(
\mathrm{ip}^{(2,X)@}_{s}\left(
P''
\right)\right)&=
\mathrm{sc}^{(0,2)}_{s}\left(
\mathrm{ip}^{(2,X)@}_{s}\left(
Q''
\right)\right)
\\
\mathrm{tg}^{(0,2)}_{s}\left(
\mathrm{ip}^{(2,X)@}_{s}\left(
P'
\right)\right)&=
\mathrm{tg}^{(0,2)}_{s}\left(
\mathrm{ip}^{(2,X)@}_{s}\left(
Q'
\right)\right);
\end{align*}
finally the third equivalence follows from left to right because $\circ^{0}_{s}$ is a defined operation in $\mathbf{Pth}_{\boldsymbol{\mathcal{A}}^{(2)}}$, according to Proposition~\ref{PDPthCatAlg}, and $\mathrm{ip}^{(2,X)@}$ is a $\Sigma^{\boldsymbol{\mathcal{A}}^{(2)}}$-homomorphism, according to Definition~\ref{DDIp}, and from right to left by Lemma~\ref{LDThetaCongSub} and by the description of the $0$-composition operation in $\mathbf{Pth}_{\boldsymbol{\mathcal{A}}^{(2)}}$, according to Proposition~\ref{PDPthCatAlg}.

Assume, without loss of generality, that $\mathrm{ip}^{(2,X)@}_{s}(P)$ is a second-order path in $\mathrm{Pth}_{\boldsymbol{\mathcal{A}}^{(2)},s}$. As we have seen before, this is the case exactly when $\mathrm{ip}^{(2,X)@}_{s}(P')$ and $\mathrm{ip}^{(2,X)@}_{s}(P'')$ are second-order paths in $\mathrm{Pth}_{\boldsymbol{\mathcal{A}}^{(2)},s}$ satisfying that 
\[
\mathrm{sc}^{(0,2)}_{s}\left(
\mathrm{ip}^{(2,X)@}_{s}\left(
P''
\right)\right)=
\mathrm{tg}^{(0,2)}_{s}\left(
\mathrm{ip}^{(2,X)@}_{s}\left(
P'
\right)\right).
\]
By induction, we also have that $\mathrm{ip}^{(2,X)@}_{s}(Q')$ and $\mathrm{ip}^{(2,X)@}_{s}(Q'')$ are second-order paths in $\mathrm{Pth}_{\boldsymbol{\mathcal{A}}^{(2)},s}$. As we have proven before, their $0$-composition is well-defined. Moreover, by induction, the following equality holds
\allowdisplaybreaks
\begin{align*}
\left(
\mathrm{ip}^{(2,X)@}_{s}\left(
P'
\right)
,
\mathrm{ip}^{(2,X)@}_{s}\left(
Q'
\right)
\right)&\in\Upsilon^{[1]}_{s};
\\
\left(
\mathrm{ip}^{(2,X)@}_{s}\left(
P''
\right)
,
\mathrm{ip}^{(2,X)@}_{s}\left(
Q''
\right)
\right)&\in\Upsilon^{[1]}_{s}.
\end{align*}

Note that the following chain of equalities holds
\allowdisplaybreaks
\begin{align*}
\mathrm{pr}^{\Upsilon^{[1]}}_{s}\left(
\mathrm{ip}^{(2,X)@}_{s}\left(
P
\right)
\right)
&=
\mathrm{pr}^{\Upsilon^{[1]}}_{s}\left(
\mathrm{ip}^{(2,X)@}_{s}\left(
P''
\circ^{0\mathbf{T}_{\Sigma^{\boldsymbol{\mathcal{A}}^{(2)}}}(X)}_{s}
P'
\right)
\right)
\tag{1}
\\&=
\mathrm{pr}^{\Upsilon^{[1]}}_{s}\left(
\mathrm{ip}^{(2,X)@}_{s}\left(
P''
\right)
\circ^{0\mathbf{Pth}_{\boldsymbol{\mathcal{A}}^{(2)}}}_{s}
\mathrm{ip}^{(2,X)@}_{s}\left(
P'
\right)
\right)
\tag{2}
\\&=
\mathrm{pr}^{\Upsilon^{[1]}}_{s}\left(
\mathrm{ip}^{(2,X)@}_{s}\left(
Q''
\right)
\circ^{0\mathbf{Pth}_{\boldsymbol{\mathcal{A}}^{(2)}}}_{s}
\mathrm{ip}^{(2,X)@}_{s}\left(
Q'
\right)
\right)
\tag{3}
\\&=
\mathrm{pr}^{\Upsilon^{[1]}}_{s}\left(
\mathrm{ip}^{(2,X)@}_{s}\left(
Q''
\circ^{0\mathbf{T}_{\Sigma^{\boldsymbol{\mathcal{A}}^{(2)}}}(X)}_{s}
Q'
\right)
\right)
\tag{4}
\\&=
\mathrm{pr}^{\Upsilon^{[1]}}_{s}\left(
\mathrm{ip}^{(2,X)@}_{s}\left(
Q
\right)
\right).
\tag{5}
\end{align*}

In the just stated chain of equalities, the first equality unravels the description of $P$; the second equivalence follows from the fact that $\circ^{0}_{s}$ is a defined operation in $\mathbf{Pth}_{\boldsymbol{\mathcal{A}}^{(2)}}$, according to Proposition~\ref{PDPthCatAlg},  and $\mathrm{ip}^{(2,X)@}$ is a $\Sigma^{\boldsymbol{\mathcal{A}}^{(2)}}$-homomorphism, according to Definition~\ref{DDIp}, and the fact that we are assuming that $\mathrm{ip}^{(2,X)@}_{s}(P)$ is a second-order path in $\mathrm{Pth}_{\boldsymbol{\mathcal{A}}^{(2)},s}$; the third equality follows by induction and from the fact that, according to Definition~\ref{DDUpsCong}, $\Upsilon^{[1]}$ is a $\Sigma^{\boldsymbol{\mathcal{A}}^{(2)}}$-congruence; the fourth equality follows from the fact that $\circ^{0}_{s}$ is a defined operation in $\mathbf{Pth}_{\boldsymbol{\mathcal{A}}^{(2)}}$, according to Proposition~\ref{PDPthCatAlg},  and $\mathrm{ip}^{(2,X)@}$ is a $\Sigma^{\boldsymbol{\mathcal{A}}^{(2)}}$-homomorphism, according to Definition~\ref{DDIp}, and the fact that we are assuming that $\mathrm{ip}^{(2,X)@}_{s}(Q')$ and $\mathrm{ip}^{(2,X)@}_{s}(Q'')$ are second-order paths in $\mathrm{Pth}_{\boldsymbol{\mathcal{A}}^{(2)},s}$ whose $0$-composition is well-defined; finally, the last equality recovers the description of $Q$.

This completes the case of the $0$-composition.

\textsf{$\gamma$ is the $1$-source operation symbol.}

Let $\gamma$ be the operation symbol $\mathrm{sc}^{1}_{s}$ in $\Sigma^{\boldsymbol{\mathcal{A}}^{(2)}}_{s,s}$ and let $(P',Q')$ be a family of pairs in $\mathrm{C}^{n}_{\Sigma^{\boldsymbol{\mathcal{A}}^{(2)}}}(\Psi^{(1)})_{s}$ for which 
$$
(P,Q)=
\left(
\mathrm{sc}^{1\mathbf{T}_{\Sigma^{\boldsymbol{\mathcal{A}}^{(2)}}}(X)}_{s}\left(
P'
\right)
,
\mathrm{sc}^{1\mathbf{T}_{\Sigma^{\boldsymbol{\mathcal{A}}^{(2)}}}(X)}_{s}\left(
Q'
\right)
\right).
$$

Note that the following chain of equivalences holds
\allowdisplaybreaks
\begin{align*}
\mathrm{ip}^{(2,X)@}_{s}(P)\in\mathrm{Pth}_{\boldsymbol{\mathcal{A}}^{(2)},s}
&\Leftrightarrow
\mathrm{ip}^{(2,X)@}_{s}(P')\in\mathrm{Pth}_{\boldsymbol{\mathcal{A}}^{(2)},s}
\tag{1}
\\&\Leftrightarrow
\mathrm{ip}^{(2,X)@}_{s}(Q')\in\mathrm{Pth}_{\boldsymbol{\mathcal{A}}^{(2)},s}
\tag{2}
\\&\Leftrightarrow
\mathrm{ip}^{(2,X)@}_{s}(Q)\in\mathrm{Pth}_{\boldsymbol{\mathcal{A}}^{(2)},s}.
\tag{3}
\end{align*}

In the just stated chain of equivalences, the first equivalence follows from left to right by Lemma~\ref{LDThetaCongSub} and from right to left because $\mathrm{sc}^{1}_{s}$ is a total operation in $\mathbf{Pth}_{\boldsymbol{\mathcal{A}}^{(2)}}$, according to Proposition~\ref{PDPthDCatAlg}, and $\mathrm{ip}^{(2,X)@}$ is a $\Sigma^{\boldsymbol{\mathcal{A}}^{(2)}}$-homomorphism, according to Definition~\ref{DDIp}; the second equivalence follows by induction, finally the third equivalence follows from left to right because $\mathrm{sc}^{1}_{s}$ is a total operation in $\mathbf{Pth}_{\boldsymbol{\mathcal{A}}^{(2)}}$, according to Proposition~\ref{PDPthAlg},  and $\mathrm{ip}^{(2,X)@}$ is a $\Sigma^{\boldsymbol{\mathcal{A}}^{(2)}}$-homomorphism, according to Definition~\ref{DDIp}, and from right to left by Lemma~\ref{LDThetaCongSub}.

Assume, without loss of generality, that $\mathrm{ip}^{(2,X)@}_{s}(P)$ is a second-order path in $\mathrm{Pth}_{\boldsymbol{\mathcal{A}}^{(2)},s}$. As we have seen before, this is the case exactly when $\mathrm{ip}^{(2,X)@}_{s}(P')$ is a second-order path in $\mathrm{Pth}_{\boldsymbol{\mathcal{A}}^{(2)},s}$. By induction, we also have that $\mathrm{ip}^{(2,X)@}_{s}(Q')$ is a second-order path in $\mathrm{Pth}_{\boldsymbol{\mathcal{A}}^{(2)},s}$. Moreover, by induction, the following equality holds
\[
\left(
\mathrm{ip}^{(2,X)@}_{s}\left(
P'
\right),
\mathrm{ip}^{(2,X)@}_{s}\left(
Q'
\right)\right)\in\Upsilon^{[1]}_{s}.
\]

Note that the following chain of equalities holds
\allowdisplaybreaks
\begin{align*}
\mathrm{pr}^{\Upsilon^{[1]}}_{s}\left(
\mathrm{ip}^{(2,X)@}_{s}\left(
P
\right)
\right)
&=
\mathrm{pr}^{\Upsilon^{[1]}}_{s}\left(
\mathrm{ip}^{(2,X)@}_{s}\left(
\mathrm{sc}^{1\mathbf{T}_{\Sigma^{\boldsymbol{\mathcal{A}}^{(2)}}}(X)}_{s}
\left(
P'
\right)
\right)
\right)
\tag{1}
\\&=
\mathrm{pr}^{\Upsilon^{[1]}}_{s}\left(
\mathrm{sc}^{1\mathbf{Pth}_{\boldsymbol{\mathcal{A}}^{(2)}}}_{s}
\left(
\mathrm{ip}^{(2,X)@}_{s}\left(
P'
\right)
\right)
\right)
\tag{2}
\\&=
\mathrm{pr}^{\Upsilon^{[1]}}_{s}\left(
\mathrm{sc}^{1\mathbf{Pth}_{\boldsymbol{\mathcal{A}}^{(2)}}}_{s}
\left(
\mathrm{ip}^{(2,X)@}_{s}\left(
Q'
\right)
\right)
\right)
\tag{3}
\\&=
\mathrm{pr}^{\Upsilon^{[1]}}_{s}\left(
\mathrm{ip}^{(2,X)@}_{s}\left(
\mathrm{sc}^{1\mathbf{T}_{\Sigma^{\boldsymbol{\mathcal{A}}^{(2)}}}(X)}_{s}
\left(
Q'
\right)
\right)
\right)
\tag{4}
\\&=
\mathrm{pr}^{\Upsilon^{[1]}}_{s}\left(
\mathrm{ip}^{(2,X)@}_{s}\left(
Q
\right)
\right).
\tag{5}
\end{align*}

In the just stated chain of equalities, the first equality unravels the description of $P$; the second equivalence follows from the fact that $\mathrm{sc}^{1}_{s}$ is a total operation in $\mathbf{Pth}_{\boldsymbol{\mathcal{A}}^{(2)}}$, according to Proposition~\ref{PDPthDCatAlg},  and $\mathrm{ip}^{(2,X)@}$ is a $\Sigma^{\boldsymbol{\mathcal{A}}^{(2)}}$-homomorphism, according to Definition~\ref{DDIp}, and the fact that we are assuming that $\mathrm{ip}^{(2,X)@}_{s}(P)$ is a second-order path in $\mathrm{Pth}_{\boldsymbol{\mathcal{A}}^{(2)},s}$; the third equality follows by induction and from the fact that, according to Definition~\ref{DDUpsCong}, $\Upsilon^{[1]}$ is a $\Sigma^{\boldsymbol{\mathcal{A}}^{(2)}}$-congruence; the fourth equality follows from the fact that $\mathrm{sc}^{1}_{s}$ is a total operation in $\mathbf{Pth}_{\boldsymbol{\mathcal{A}}^{(2)}}$, according to Proposition~\ref{PDPthDCatAlg},  and $\mathrm{ip}^{(2,X)@}$ is a $\Sigma^{\boldsymbol{\mathcal{A}}^{(2)}}$-homomorphism, according to Definition~\ref{DDIp}, and the fact that we are assuming that $\mathrm{ip}^{(2,X)@}_{s}(Q')$ is a second-order path in $\mathrm{Pth}_{\boldsymbol{\mathcal{A}}^{(2)},s}$; finally, the last equality recovers the description of $Q$.

This completes the case of the $1$-source.

\textsf{$\gamma$ is the $1$-target operation symbol.}

Let $\gamma$ be the operation symbol $\mathrm{tg}^{1}_{s}$ in $\Sigma^{\boldsymbol{\mathcal{A}}^{(2)}}_{s,s}$ and let $(P',Q')$ be a family of pairs in $\mathrm{C}^{n}_{\Sigma^{\boldsymbol{\mathcal{A}}^{(2)}}}(\Psi^{(1)})_{s}$ for which 
$$
(P,Q)=
\left(
\mathrm{tg}^{1\mathbf{T}_{\Sigma^{\boldsymbol{\mathcal{A}}^{(2)}}}(X)}_{s}\left(
P'
\right)
,
\mathrm{tg}^{1\mathbf{T}_{\Sigma^{\boldsymbol{\mathcal{A}}^{(2)}}}(X)}_{s}\left(
Q'
\right)
\right).
$$

Then, the following properties hold
\begin{itemize}
\item[(i)] $\mathrm{ip}^{(2,X)@}_{s}(P)\in\mathrm{Pth}_{\boldsymbol{\mathcal{A}}^{(2)},s}$ if, and only if, $\mathrm{ip}^{(2,X)@}_{s}(Q)\in\mathrm{Pth}_{\boldsymbol{\mathcal{A}}^{(2)},s}$;
\item[(ii)] If $\mathrm{ip}^{(2,X)@}_{s}(P)\in\mathrm{Pth}_{\boldsymbol{\mathcal{A}}^{(2)},s}$ or  $\mathrm{ip}^{(2,X)@}_{s}(Q)\in\mathrm{Pth}_{\boldsymbol{\mathcal{A}}^{(2)},s}$ is a second-order path in $\mathrm{Pth}_{\boldsymbol{\mathcal{A}}^{(2)},s}$ then
$
(
\mathrm{ip}^{(2,X)@}_{s}(
P),
\mathrm{ip}^{(2,X)@}_{s}(
Q))\in\Upsilon^{[1]}_{s}.
$
\end{itemize}

The proof of this case is similar to that of the $1$-source.

This completes the case of the $1$-target.

\textsf{$\gamma$ is the $1$-composition operation symbol.}

Let $\gamma$ be the operation symbol $\circ^{1}_{s}$ in $\Sigma^{\boldsymbol{\mathcal{A}}^{(2)}}_{s,s}$ and let $(P',Q')$ and $(P'',Q'')$ be two families of pairs in $\mathrm{C}^{n}_{\Sigma^{\boldsymbol{\mathcal{A}}^{(2)}}}(\Psi^{(1)})_{s}$ for which 
$$
(P,Q)=
\left(
P''
\circ^{1\mathbf{T}_{\Sigma^{\boldsymbol{\mathcal{A}}^{(2)}}}(X)}_{s}
P'
,
Q''
\circ^{1\mathbf{T}_{\Sigma^{\boldsymbol{\mathcal{A}}^{(2)}}}(X)}_{s}
Q'
\right).
$$

Note that the following chain of equivalences holds
\allowdisplaybreaks
\begin{align*}
\mathrm{ip}^{(2,X)@}_{s}(P)\in\mathrm{Pth}_{\boldsymbol{\mathcal{A}}^{(2)},s}
&\Leftrightarrow
\left\lbrace
\begin{array}{l}
\mathrm{ip}^{(2,X)@}_{s}(P'),\mathrm{ip}^{(2,X)@}_{s}(P'')\in\mathrm{Pth}_{\boldsymbol{\mathcal{A}}^{(2)},s}
\\
\mathrm{sc}^{([1],2)}_{s}(
\mathrm{ip}^{(2,X)@}_{s}(P'')
)
=
\mathrm{tg}^{([1],2)}_{s}(
\mathrm{ip}^{(2,X)@}_{s}(P')
)
\end{array}
\right.
\tag{1}
\\&\Leftrightarrow
\left\lbrace
\begin{array}{l}
\mathrm{ip}^{(2,X)@}_{s}(Q'),\mathrm{ip}^{(2,X)@}_{s}(Q'')\in\mathrm{Pth}_{\boldsymbol{\mathcal{A}}^{(2)},s}
\\
\mathrm{sc}^{([1],2)}_{s}(
\mathrm{ip}^{(2,X)@}_{s}(Q'')
)
=
\mathrm{tg}^{([1],2)}_{s}(
\mathrm{ip}^{(2,X)@}_{s}(Q')
)
\end{array}
\right.
\tag{2}
\\&\Leftrightarrow
\mathrm{ip}^{(2,X)@}_{s}(Q)\in\mathrm{Pth}_{\boldsymbol{\mathcal{A}}^{(2)},s}.
\tag{3}
\end{align*}

In the just stated chain of equivalences, the first equivalence follows from left to right by Lemma~\ref{LDThetaCongSub} and by the description of the $1$-composition operation in $\mathbf{Pth}_{\boldsymbol{\mathcal{A}}^{(2)}}$, according to Proposition~\ref{PDPthDCatAlg}, and from right to left because $\circ^{1}_{s}$ is a defined operation in $\mathbf{Pth}_{\boldsymbol{\mathcal{A}}^{(2)}}$, according to Proposition~\ref{PDPthDCatAlg}, and $\mathrm{ip}^{(2,X)@}$ is a $\Sigma^{\boldsymbol{\mathcal{A}}^{(2)}}$-homomorphism, according to Definition~\ref{DDIp}; the second equivalence follows by induction. Note also that
\allowdisplaybreaks
\begin{align*}
\left(
\mathrm{ip}^{(2,X)@}_{s}\left(
P'
\right),
\mathrm{ip}^{(2,X)@}_{s}\left(
Q'
\right)
\right)&\in\Upsilon^{[1]}_{s};
\\
\left(
\mathrm{ip}^{(2,X)@}_{s}\left(
P''
\right),
\mathrm{ip}^{(2,X)@}_{s}\left(
Q''
\right)
\right)&\in\Upsilon^{[1]}_{s}.
\end{align*}
Therefore, according to Lemma~\ref{LDUpsCong}, we have that 
\allowdisplaybreaks
\begin{align*}
\mathrm{sc}^{([1],2)}_{s}\left(
\mathrm{ip}^{(2,X)@}_{s}\left(
P''
\right)\right)&=
\mathrm{sc}^{([1],2)}_{s}\left(
\mathrm{ip}^{(2,X)@}_{s}\left(
Q''
\right)\right)
\\
\mathrm{tg}^{([1],2)}_{s}\left(
\mathrm{ip}^{(2,X)@}_{s}\left(
P'
\right)\right)&=
\mathrm{tg}^{([1],2)}_{s}\left(
\mathrm{ip}^{(2,X)@}_{s}\left(
Q'
\right)\right);
\end{align*}
finally the third equivalence follows from left to right because $\circ^{1}_{s}$ is a defined operation in $\mathbf{Pth}_{\boldsymbol{\mathcal{A}}^{(2)}}$, according to Proposition~\ref{PDPthDCatAlg}, and $\mathrm{ip}^{(2,X)@}$ is a $\Sigma^{\boldsymbol{\mathcal{A}}^{(2)}}$-homomorphism, according to Definition~\ref{DDIp}, and from right to left by Lemma~\ref{LDThetaCongSub} and by the description of the $1$-composition operation in $\mathbf{Pth}_{\boldsymbol{\mathcal{A}}^{(2)}}$, according to Proposition~\ref{PDPthDCatAlg}.

Assume, without loss of generality, that $\mathrm{ip}^{(2,X)@}_{s}(P)$ is a second-order path in $\mathrm{Pth}_{\boldsymbol{\mathcal{A}}^{(2)},s}$. As we have seen before, this is the case exactly when $\mathrm{ip}^{(2,X)@}_{s}(P')$ and $\mathrm{ip}^{(2,X)@}_{s}(P'')$ are second-order paths in $\mathrm{Pth}_{\boldsymbol{\mathcal{A}}^{(2)},s}$ satisfying that 
\[
\mathrm{sc}^{([1],2)}_{s}\left(
\mathrm{ip}^{(2,X)@}_{s}\left(
P''
\right)\right)=
\mathrm{tg}^{([1],2)}_{s}\left(
\mathrm{ip}^{(2,X)@}_{s}\left(
P'
\right)\right).
\]
By induction, we also have that $\mathrm{ip}^{(2,X)@}_{s}(Q')$ and $\mathrm{ip}^{(2,X)@}_{s}(Q'')$ are second-order paths in $\mathrm{Pth}_{\boldsymbol{\mathcal{A}}^{(2)},s}$. As we have proven before, their $1$-composition is well-defined. Moreover, by induction, the following equality holds
\allowdisplaybreaks
\begin{align*}
\left(
\mathrm{ip}^{(2,X)@}_{s}\left(
P'
\right),
\mathrm{ip}^{(2,X)@}_{s}\left(
Q'
\right)\right)&\in\Upsilon^{[1]}_{s}.
\\
\left(
\mathrm{ip}^{(2,X)@}_{s}\left(
P''
\right),
\mathrm{ip}^{(2,X)@}_{s}\left(
Q''
\right)\right)&\in\Upsilon^{[1]}_{s}.
\end{align*}

Note that the following chain of equalities holds
\allowdisplaybreaks
\begin{align*}
\mathrm{pr}^{\Upsilon^{[1]}}_{s}\left(
\mathrm{ip}^{(2,X)@}_{s}\left(
P
\right)
\right)
&=
\mathrm{pr}^{\Upsilon^{[1]}}_{s}\left(
\mathrm{ip}^{(2,X)@}_{s}\left(
P''
\circ^{1\mathbf{T}_{\Sigma^{\boldsymbol{\mathcal{A}}^{(2)}}}(X)}_{s}
P'
\right)
\right)
\tag{1}
\\&=
\mathrm{pr}^{\Upsilon^{[1]}}_{s}\left(
\mathrm{ip}^{(2,X)@}_{s}\left(
P''
\right)
\circ^{1\mathbf{Pth}_{\boldsymbol{\mathcal{A}}^{(2)}}}_{s}
\mathrm{ip}^{(2,X)@}_{s}\left(
P'
\right)
\right)
\tag{2}
\\&=
\mathrm{pr}^{\Upsilon^{[1]}}_{s}\left(
\mathrm{ip}^{(2,X)@}_{s}\left(
Q''
\right)
\circ^{1\mathbf{Pth}_{\boldsymbol{\mathcal{A}}^{(2)}}}_{s}
\mathrm{ip}^{(2,X)@}_{s}\left(
Q'
\right)
\right)
\tag{3}
\\&=
\mathrm{pr}^{\Upsilon^{[1]}}_{s}\left(
\mathrm{ip}^{(2,X)@}_{s}\left(
Q''
\circ^{1\mathbf{T}_{\Sigma^{\boldsymbol{\mathcal{A}}^{(2)}}}(X)}_{s}
Q'
\right)
\right)
\tag{4}
\\&=
\mathrm{pr}^{\Upsilon^{[1]}}_{s}\left(
\mathrm{ip}^{(2,X)@}_{s}\left(
Q
\right)
\right).
\tag{5}
\end{align*}

In the just stated chain of equalities, the first equality unravels the description of $P$; the second equivalence follows from the fact that $\circ^{1}_{s}$ is a defined operation in $\mathbf{Pth}_{\boldsymbol{\mathcal{A}}^{(2)}}$, according to Proposition~\ref{PDPthDCatAlg},  and $\mathrm{ip}^{(2,X)@}$ is a $\Sigma^{\boldsymbol{\mathcal{A}}^{(2)}}$-homomorphism, according to Definition~\ref{DDIp}, and the fact that we are assuming that $\mathrm{ip}^{(2,X)@}_{s}(P)$ is a second-order path in $\mathrm{Pth}_{\boldsymbol{\mathcal{A}}^{(2)},s}$; the third equality follows by induction and from the fact that, according to Definition~\ref{DDUpsCong}, $\Upsilon^{[1]}$ is a $\Sigma^{\boldsymbol{\mathcal{A}}^{(2)}}$-congruence; the fourth equality follows from the fact that $\circ^{1}_{s}$ is a defined operation in $\mathbf{Pth}_{\boldsymbol{\mathcal{A}}^{(2)}}$, according to Proposition~\ref{PDPthCatAlg},  and $\mathrm{ip}^{(2,X)@}$ is a $\Sigma^{\boldsymbol{\mathcal{A}}^{(2)}}$-homomorphism, according to Definition~\ref{DDIp}, and the fact that we are assuming that $\mathrm{ip}^{(2,X)@}_{s}(Q')$ and $\mathrm{ip}^{(2,X)@}_{s}(Q'')$ are second-order paths in $\mathrm{Pth}_{\boldsymbol{\mathcal{A}}^{(2)},s}$ whose $1$-composition is well-defined; finally, the last equality recovers the description of $Q$.

This completes the case of the $1$-composition.

This completes the proof.
\end{proof}

\chapter{
\texorpdfstring
{The congruence $\Theta^{[2]}\vee \Psi^{[1]}$ on $\mathbf{T}_{\Sigma^{\boldsymbol{\mathcal{A}}^{(2)}}}(X)$ and its relatives}
{The second-order congruence on terms}
}\label{S2L}

In this chapter we consider the congruence $\Theta^{[2]}\vee \Psi^{[1]}$ on $\mathbf{T}_{\Sigma^{\boldsymbol{\mathcal{A}}^{(2)}}}(X)$, i.e., the supremum of the $\Sigma^{\boldsymbol{\mathcal{A}}^{(2)}}$-congruences $\Theta^{[2]}$ and $\Psi^{[1]}$ on $\mathbf{T}_{\Sigma^{\boldsymbol{\mathcal{A}}^{(2)}}}(X)$, that we will denote by $\Theta^{\llbracket 2 \rrbracket }$. We will also denote by $\llbracket \mathrm{T}_{\Sigma^{\boldsymbol{\mathcal{A}}^{(2)}}}(X)\rrbracket_{\Theta^{\llbracket 2\rrbracket}}$ the corresponding quotient and the $\Theta^{\llbracket 2 \rrbracket }$-equivalence class of a term $P$ in $\mathrm{T}_{\Sigma^{\boldsymbol{\mathcal{A}}^{(2)}}}(X)_{s}$ will be denoted by $\llbracket P\rrbracket_{\Theta^{\llbracket 2 \rrbracket}_{s}}$. After this we prove that if two second-order paths are $\Upsilon^{[1]}_{s}$-related then their images under the second-order Curry-Howard mapping are $\Theta^{\llbracket 2 \rrbracket }_{s}$-related. Moreover, if two second-order paths are $\mathrm{Ker}(\mathrm{CH}^{(2)})_{s}\vee \Upsilon^{[1]}_{s}$-related then their images under the second-order Curry-Howard mapping are $\Theta^{\llbracket 2 \rrbracket }_{s}$-related. We also prove that if a term $P$ in $\mathrm{T}_{\Sigma^{\boldsymbol{\mathcal{A}}^{(2)}}}(X)_{s}$ is such that $\mathrm{ip}^{(2,X)@}_{s}(P)$ is a second-order path in $\mathrm{Pth}_{\boldsymbol{\mathcal{A}}^{(2)},s}$, then the terms $P$ and $\mathrm{CH}^{(2)}_{s}(\mathrm{ip}^{(2,X)@}_{s}(P))$ are $\Theta^{\llbracket 2 \rrbracket}_{s}$-related. Moreover, if a pair of terms in $\mathrm{T}_{\Sigma^{\boldsymbol{\mathcal{A}}^{(2)}}}(X)_{s}$ are $\Theta^{\llbracket 2 \rrbracket }_{s}$-related and, after applying the $\mathrm{ip}^{(2,X)@}$ mapping, they become second-order paths, then these are $\mathrm{Ker}(\mathrm{CH}^{(2)})_{s}\vee \Upsilon^{[1]}_{s}$-related. This, in particular implies that, if a term $P$ is $\Theta^{\llbracket 2\rrbracket}_{s}$-related with $\mathrm{CH}^{(2)}_{s}(\mathfrak{P}^{(2)})$, then $\mathrm{ip}^{(2,X)@}_{s}(P)$  is itself a second-order path in the equivalence class $\llbracket \mathfrak{P}^{(2)}\rrbracket_{s}$. Furthermore it is shown that, for two second-order paths  $\mathfrak{P}^{(2)}$ and $\mathfrak{P}'^{(2)}$ in $\mathrm{Pth}_{\boldsymbol{\mathcal{A}}^{(2)},s}$, if $\mathrm{CH}^{(2)}_{s}(\mathfrak{P}^{(2)})$ and $\mathrm{CH}^{(2)}_{s}(\mathfrak{P}'^{(2)})$ are  $\Theta^{\llbracket 2 \rrbracket}_{s}$-related then the equivalence classes $\llbracket \mathfrak{P}^{(2)}\rrbracket_{s}$ and  $\llbracket\mathfrak{P}'^{(2)}\rrbracket_{s}$ are equal.


\begin{restatable}{definition}{DDVCong}
\label{DDVCong}
\index{Theta!second-order!$\Theta^{\llbracket 2 \rrbracket}$}
We define  $\Theta^{\llbracket 2 \rrbracket}$ to be the relation $\Theta^{[2]}\vee \Psi^{[1]}$ in $\mathrm{T}_{\Sigma^{\boldsymbol{\mathcal{A}}^{(2)}}}(X)$.
\end{restatable}

 Let us recall that, for a sort $s\in S$, two terms $P$ and $Q$ in $\mathrm{T}_{\Sigma^{\boldsymbol{\mathcal{A}}}}(X)_{s}$ are related under $\Theta^{\llbracket 2\rrbracket}_{s}$ if, and only if, there exists a natural number $m\in \mathbb{N}-\{0\}$ and a sequence of terms $(R_{k})_{k\in m+1}$ in $\mathrm{T}_{\Sigma^{\boldsymbol{\mathcal{A}}^{(2)}}}(X)_{s}^{m+1}$ satisfying that 
\begin{enumerate}
\item $R_{0}=P$;
\item $R_{m}=Q$;
\item for every $k\in m$, 
$(R_{k},R_{k+1})\in 
 \Theta^{[2]}_{s}\cup \Psi^{[1]}_{s}$.
\end{enumerate}

\begin{restatable}{convention}{CDVTermConv}
\label{CDVTermConv} 
\index{terms!second-order!$\llbracket   P\rrbracket_{\Theta^{\llbracket 2\rrbracket}_{s}}$}
\index{terms!second-order!$\llbracket \mathrm{T}_{\Sigma^{\boldsymbol{\mathcal{A}}^{(2)}}}(X)\rrbracket_{\Theta^{\llbracket 2\rrbracket}}$}
In order to simplify the presentation, for a sort $s\in S$ and a term $P\in\mathrm{T}_{\Sigma^{\boldsymbol{\mathcal{A}}^{(2)}}}(X)_{s}$, we will let $\llbracket   P\rrbracket_{\Theta^{\llbracket 2\rrbracket}_{s}}$ stand for $[P]_{\Theta^{[2]}_{s}\vee \Psi^{[1]}_{s}}$, the $\Theta^{[2]}_{s}\vee \Psi^{[1]}_{s}$-class of $P$, and we will denote by $\llbracket \mathrm{T}_{\Sigma^{\boldsymbol{\mathcal{A}}^{(2)}}}(X)\rrbracket_{\Theta^{\llbracket 2\rrbracket}}$ the many-sorted quotient
\[
\llbracket \mathrm{T}_{\Sigma^{\boldsymbol{\mathcal{A}}^{(2)}}}(X)\rrbracket_{\Theta^{\llbracket 2\rrbracket}}=
\mathrm{T}_{\Sigma^{\boldsymbol{\mathcal{A}}^{(2)}}}(X)/{\Theta^{[2]}\vee \Psi^{[1]}}.
\]
\end{restatable}

In the following proposition we will prove that, for every sort $s\in S$, if two second-order paths in $\mathrm{Pth}_{\boldsymbol{\mathcal{A}}^{(2)},s}$ are $\Upsilon^{[1]}_{s}$-related, then their respective values under the second-order Curry-Howard mapping are $\Theta^{\llbracket 2\rrbracket}_{s}$-related.

\begin{restatable}{proposition}{PDUpsDCH}
\label{PDUpsDCH} Let $s$ be a sort in $S$ and let $\mathfrak{P}^{(2)}$ and $\mathfrak{Q}^{(2)}$ be second-order paths in $\mathrm{Pth}_{\boldsymbol{\mathcal{A}}^{(2)},s}$ satisfying that $[\mathfrak{P}^{(2)}]_{\Upsilon^{[1]}_{s}}=[\mathfrak{Q}^{(2)}]_{\Upsilon^{[1]}_{s}}$, then
\[
\Bigl\llbracket
\mathrm{CH}^{(2)}_{s}\left(
\mathfrak{P}^{(2)}
\right)
\Bigr\rrbracket_{\Theta^{\llbracket 2\rrbracket}_{s}}
=
\Bigl\llbracket
\mathrm{CH}^{(2)}_{s}\left(
\mathfrak{Q}^{(2)}
\right)
\Bigr\rrbracket_{\Theta^{\llbracket 2\rrbracket}_{s}}.
\]
\end{restatable}
\begin{proof}
We recall from Remark~\ref{RDUpsCong} that 
\[
\Upsilon^{[1]}=\bigcup_{n\in\mathbb{N}} \mathrm{C}^{n}_{\mathbf{Pth}_{\boldsymbol{\mathcal{A}}^{(2)}}}\left(\Upsilon^{(1)}\right).
\]

We prove the statement by induction on $n\in\mathbb{N}$.

\textsf{Base step of the induction.}

Let us recall from Remark~\ref{RDUpsCong} that 
\[
\mathrm{C}^{0}_{\mathbf{Pth}_{\boldsymbol{\mathcal{A}}^{(2)}}}\left(\Upsilon^{(1)}\right)
=
\left(
\Upsilon^{(1)}
\right)
\cup
\left(
\Upsilon^{(1)}
\right)^{-1}
\cup
\Delta_{\mathrm{Pth}_{\boldsymbol{\mathcal{A}}^{(2)}}}.
\]

The statement trivially holds for a pair $(\mathfrak{P}^{(2)},\mathfrak{Q}^{(2)})$ in $\Delta_{\mathrm{Pth}_{\boldsymbol{\mathcal{A}}^{(2)}},s}$.

If the pair $(\mathfrak{P}^{(2)},\mathfrak{Q}^{(2)})$ is in $\Upsilon^{(1)}_{s}$, following Definition~\ref{DDUpsCong}, we conclude that 
\[
\left(
\mathrm{CH}^{(2)}_{s}\left(\mathfrak{P}^{(2)}\right),
\mathrm{CH}^{(2)}_{s}\left(\mathfrak{Q}^{(2)}\right)
\right)
\in \Psi^{(1)}_{s}\subseteq \Psi^{[1]}_{s}\subseteq \Theta^{\llbracket 2\rrbracket}_{s}.
\]

If the pair $(\mathfrak{Q}^{(2)},\mathfrak{P}^{(2)})$ is in $(\Upsilon^{(1)})^{-1}_{s}$, following Definition~\ref{DDUpsCong}, we conclude that 
\[
\left(
\mathrm{CH}^{(2)}_{s}\left(\mathfrak{Q}^{(2)}\right),
\mathrm{CH}^{(2)}_{s}\left(\mathfrak{P}^{(2)}\right)
\right)
\in \left(\Psi^{(1)}\right)^{-1}_{s}\subseteq \Psi^{[1]}_{s}\subseteq \Theta^{\llbracket 2\rrbracket}_{s}.
\]

This completes the base case.

\textsf{Inductive step of the induction.}

Assume the statement holds for $n\in\mathbb{N}$, i.e., for every sort $s\in S$ and every pair of second-order paths $(\mathfrak{P}^{(2)},\mathfrak{Q}^{(2)})$ in $\mathrm{C}^{n}_{\mathrm{Pth}_{\boldsymbol{\mathcal{A}}^{(2)}}}(\Upsilon	^{(1)})_{s}$ then
\[
\Bigl\llbracket
\mathrm{CH}^{(2)}_{s}\left(
\mathfrak{P}^{(2)}
\right)
\Bigr\rrbracket_{\Theta^{\llbracket 2\rrbracket}_{s}}
=
\Bigl\llbracket
\mathrm{CH}^{(2)}_{s}\left(
\mathfrak{Q}^{(2)}
\right)
\Bigr\rrbracket_{\Theta^{\llbracket 2\rrbracket}_{s}}.
\]

We prove the statement for $n+1$. Let $s$ be a sort in $S$ and let $(\mathfrak{P}^{(2)},\mathfrak{Q}^{(2)})$ be a pair of second-order paths in $\mathrm{Pth}_{\boldsymbol{\mathcal{A}}^{(2)},s}$ such that $(\mathfrak{P}^{(2)},\mathfrak{Q}^{(2)})$ is a pair in $\mathrm{C}^{n+1}_{\mathbf{Pth}_{\boldsymbol{\mathcal{A}}^{(2)}}}(\Upsilon^{(1)})_{s}$. Let us recall from  Remark~\ref{RDUpsCong} that 
\allowdisplaybreaks
\begin{multline*}
\mathrm{C}^{n+1}_{\mathbf{Pth}_{\boldsymbol{\mathcal{A}}^{(2)}}}\left(\Upsilon^{(1)}\right)_{s}
=
\left(
\mathrm{C}^{n}_{\mathbf{Pth}_{\boldsymbol{\mathcal{A}}^{(2)}}}\left(\Upsilon^{(1)}\right)_{s}
\circ 
\mathrm{C}^{n}_{\mathbf{Pth}_{\boldsymbol{\mathcal{A}}^{(2)}}}\left(\Upsilon^{(1)}\right)_{s}
\right)
\cup
\\
\left(
\bigcup_{\gamma\in\Sigma^{\boldsymbol{\mathcal{A}}^{(2)}}_{\neq \lambda,s}}
\gamma^{\mathbf{Pth}_{\boldsymbol{\mathcal{A}}^{(2)}}}
\times
\gamma^{\mathbf{Pth}_{\boldsymbol{\mathcal{A}}^{(2)}}}
\left[
\mathrm{C}^{n}_{\mathbf{Pth}_{\boldsymbol{\mathcal{A}}^{(2)}}}\left(\Upsilon^{(1)}\right)_{\mathrm{ar}(\gamma)}
\right]
\right).
\end{multline*}

Then either (1) $(\mathfrak{P}^{(2)},\mathfrak{Q}^{(2)})$ is in $\mathrm{C}^{n}_{\mathbf{Pth}_{\boldsymbol{\mathcal{A}}^{(2)}}}(\Upsilon^{(1)})_{s}
\circ 
\mathrm{C}^{n}_{\mathbf{Pth}_{\boldsymbol{\mathcal{A}}^{(2)}}}(\Upsilon^{(1)})_{s}$ or (2)  $(\mathfrak{P}^{(2)},\mathfrak{Q}^{(2)})$ is in $\gamma^{\mathbf{Pth}_{\boldsymbol{\mathcal{A}}^{(2)}}}
\times
\gamma^{\mathbf{Pth}_{\boldsymbol{\mathcal{A}}^{(2)}}}
[
\mathrm{C}^{n}_{\mathbf{Pth}_{\boldsymbol{\mathcal{A}}^{(2)}}}(\Upsilon^{(1)})_{\mathrm{ar}(\gamma)}
]$, for some operation symbol $\gamma \in \Sigma^{\boldsymbol{\mathcal{A}}^{(2)}}_{\neq \lambda, s}$.

If (1), then there exists a second-order path $\mathfrak{R}^{(2)} \in\mathrm{Pth}_{\boldsymbol{\mathcal{A}}^{(2)},s}$ for which $(\mathfrak{P}^{(2)},\mathfrak{R}^{(2)})$ and $(\mathfrak{R}^{(2)},\mathfrak{Q}^{(2)})$ belong to $\mathrm{C}^{n}_{\mathbf{Pth}_{\boldsymbol{\mathcal{A}}^{(2)}},s}(\Upsilon^{(1)})_{s}$. Then, by induction, we have  that
\allowdisplaybreaks
\begin{align*}
\Bigl\llbracket
\mathrm{CH}^{(2)}_{s}\left(
\mathfrak{P}^{(2)}
\right)
\Bigr\rrbracket_{\Theta^{\llbracket 2\rrbracket}_{s}}
&=
\Bigl\llbracket
\mathrm{CH}^{(2)}_{s}\left(
\mathfrak{R}^{(2)}
\right)
\Bigr\rrbracket_{\Theta^{\llbracket 2\rrbracket}_{s}};
\\
\Bigl\llbracket
\mathrm{CH}^{(2)}_{s}\left(
\mathfrak{R}^{(2)}
\right)
\Bigr\rrbracket_{\Theta^{\llbracket 2\rrbracket}_{s}}
&=
\Bigl\llbracket
\mathrm{CH}^{(2)}_{s}\left(
\mathfrak{Q}^{(2)}
\right)
\Bigr\rrbracket_{\Theta^{\llbracket 2\rrbracket}_{s}}.
\end{align*}

Thus, by transitivity of $\Theta^{\llbracket 2\rrbracket}_{s}$, we conclude that 
\[
\Bigl\llbracket
\mathrm{CH}^{(2)}_{s}\left(
\mathfrak{P}^{(2)}
\right)
\Bigr\rrbracket_{\Theta^{\llbracket 2\rrbracket}_{s}}
=
\Bigl\llbracket
\mathrm{CH}^{(2)}_{s}\left(
\mathfrak{Q}^{(2)}
\right)
\Bigr\rrbracket_{\Theta^{\llbracket 2\rrbracket}_{s}}.
\]

If~(2), then there exists a unique word $\mathbf{s}\in S^{\star}-\{\lambda\}$, a unique operation symbol $\gamma \in\Sigma^{\boldsymbol{\mathcal{A}}^{(2)}}_{\mathbf{s},s}$ and a unique family of pairs $((\mathfrak{P}^{(2)}_{j},\mathfrak{Q}^{(2)}_{j}))_{j\in\bb{\mathbf{s}}}$ in $\mathrm{C}^{n}_{\boldsymbol{\mathcal{A}}^{(2)}}(\Upsilon^{(1)})_{\mathbf{s}}$ and in the domain of $\gamma$ for which
\[
\left(
\mathfrak{P}^{(2)},
\mathfrak{Q}^{(2)}
\right)
=
\left(
\gamma^{\mathbf{Pth}_{\boldsymbol{\mathcal{A}}^{(2)}}}
\left(
\left(
\mathfrak{P}^{(2)}_{j}
\right)_{j\in\bb{\mathbf{s}}}
\right),
\gamma^{\mathbf{Pth}_{\boldsymbol{\mathcal{A}}^{(2)}}}
\left(
\left(
\mathfrak{Q}^{(2)}_{j}
\right)_{j\in\bb{\mathbf{s}}}
\right)
\right).
\]

We will distinguish the following cases according to the different possibilities for the operation symbol $\gamma \in \Sigma^{\boldsymbol{\mathcal{A}}^{(2)}}_{\mathbf{s},s}$. Note that either $\gamma$ is an operation symbol $\sigma \in \Sigma_{\mathbf{s},s}$,  the operation symbol of $0$-source, $\mathrm{sc}^{0}_{s}$,   the operation symbol of $0$-target, $\mathrm{tg}^{0}_{s}$,   the operation symbol of $0$-composition, $\circ^{0}_{s}$,  the operation symbol of $1$-source, $\mathrm{sc}^{1}_{s}$,   the operation symbol of $1$-target, $\mathrm{tg}^{1}_{s}$, or the operation symbol of $1$-composition, $\circ^{1}_{s}$.

\textsf{$\gamma$ is an operation symbol $\sigma \in\Sigma_{\mathbf{s},s}$.}

Let $\mathbf{s}$ be a word in $S^{\star}-\{\lambda\}$, let $\sigma$ be an operation symbol in $\Sigma_{\mathbf{s},s}$ and let $((\mathfrak{P}^{(2)}_{j},\mathfrak{Q}^{(2)}_{j}))_{j\in\bb{\mathbf{s}}}$ be a family of pairs in $\mathrm{C}^{n}_{\boldsymbol{\mathcal{A}}^{(2)}}(\Upsilon^{(1)})_{\mathbf{s}}$ for which
\[
\left(
\mathfrak{P}^{(2)},
\mathfrak{Q}^{(2)}
\right)
=
\left(
\sigma^{\mathbf{Pth}_{\boldsymbol{\mathcal{A}}^{(2)}}}
\left(
\left(
\mathfrak{P}^{(2)}_{j}
\right)_{j\in\bb{\mathbf{s}}}
\right),
\sigma^{\mathbf{Pth}_{\boldsymbol{\mathcal{A}}^{(2)}}}
\left(
\left(
\mathfrak{Q}^{(2)}_{j}
\right)_{j\in\bb{\mathbf{s}}}
\right)
\right).
\]

Note that the following chain of equalities holds
\allowdisplaybreaks
\begin{align*}
\left\llbracket \mathrm{CH}^{(2)}_{s}\left(
\mathfrak{P}^{(2)}
\right)\right\rrbracket_{s}&=
\left\llbracket \mathrm{CH}^{(2)}_{s}\left(
\sigma^{\mathbf{Pth}_{\boldsymbol{\mathcal{A}}^{(2)}}}
\left(
\left(
\mathfrak{P}^{(2)}_{j}
\right)_{j\in\bb{\mathbf{s}}}
\right)
\right)\right\rrbracket_{s}
\tag{1}
\\&=
\left\llbracket 
\sigma^{\mathbf{T}_{\Sigma^{\boldsymbol{\mathcal{A}}^{(2)}}}(X)}
\left(
\left(
\mathrm{CH}^{(2)}_{s_{j}}\left(
\mathfrak{P}^{(2)}_{j}
\right)
\right)_{j\in\bb{\mathbf{s}}}
\right)
\right\rrbracket_{s}
\tag{2}
\\&=
\left\llbracket 
\sigma^{\mathbf{T}_{\Sigma^{\boldsymbol{\mathcal{A}}^{(2)}}}(X)}
\left(
\left(
\mathrm{CH}^{(2)}_{s_{j}}\left(
\mathfrak{Q}^{(2)}_{j}
\right)
\right)_{j\in\bb{\mathbf{s}}}
\right)
\right\rrbracket_{s}
\tag{3}
\\&=
\left\llbracket \mathrm{CH}^{(2)}_{s}\left(
\sigma^{\mathbf{Pth}_{\boldsymbol{\mathcal{A}}^{(2)}}}
\left(
\left(
\mathfrak{Q}^{(2)}_{j}
\right)_{j\in\bb{\mathbf{s}}}
\right)
\right)\right\rrbracket_{s}
\tag{4}
\\&=
\left\llbracket \mathrm{CH}^{(2)}_{s}\left(
\mathfrak{Q}^{(2)}
\right)\right\rrbracket_{s}.
\tag{5}
\end{align*}

In the just stated chain of equalities, the first equality unravels the description of the second-order path $\mathfrak{P}^{(2)}$; the second equality follows from the fact that $\mathrm{CH}^{(2)}$ is a $\Sigma$-homomorphism, according to Proposition~\ref{PDCHHom}; the third equality follows from the fact that $\Theta^{\llbracket 2\rrbracket}$ is a $\Sigma^{\boldsymbol{\mathcal{A}}^{(2)}}$-congruence and by the fact that, by induction, for every $j\in\bb{\mathbf{s}}$, the following equality holds
\[
\left\llbracket
\mathrm{CH}^{(2)}_{s_{j}}\left(
\mathfrak{P}^{(2)}_{j}
\right)
\right\rrbracket_{\Theta^{\llbracket 2\rrbracket}_{s_{j}}}
=
\left\llbracket
\mathrm{CH}^{(2)}_{s_{j}}\left(
\mathfrak{Q}^{(2)}_{j}
\right)
\right\rrbracket_{\Theta^{\llbracket 2\rrbracket}_{s_{j}}};
\]
the fourth equality follows from the fact that $\mathrm{CH}^{(2)}$ is a $\Sigma$-homomorphism, according to Proposition~\ref{PDCHHom}; finally, the last equality recovers  the description of the second-order path $\mathfrak{Q}^{(2)}$.

This completes the case $\sigma\in \Sigma_{\mathbf{s},s}$.

\textsf{$\gamma$ is the $0$-source operation symbol.}

Let $(\mathfrak{P}'^{(2)}, \mathfrak{Q}'^{(2)})$ be the pair in $\mathrm{C}^{n}_{\boldsymbol{\mathcal{A}}^{(2)}}(\Upsilon^{(1)})_{s}$ for which
\[
\left(
\mathfrak{P}^{(2)},
\mathfrak{Q}^{(2)}
\right)
=
\left(
\mathrm{sc}^{0\mathbf{Pth}_{\boldsymbol{\mathcal{A}}^{(2)}}}_{s}\left(
\mathfrak{P}'^{(2)}\right),
\mathrm{sc}^{0\mathbf{Pth}_{\boldsymbol{\mathcal{A}}^{(2)}}}_{s}\left(
\mathfrak{Q}'^{(2)}\right)
\right).
\]

Note that the following chain of equalities holds
\allowdisplaybreaks
\begin{align*}
\left\llbracket \mathrm{CH}^{(2)}_{s}\left(
\mathfrak{P}^{(2)}
\right)\right\rrbracket_{s}&=
\left\llbracket \mathrm{CH}^{(2)}_{s}\left(
\mathrm{sc}^{0\mathbf{Pth}_{\boldsymbol{\mathcal{A}}^{(2)}}}_{s}\left(
\mathfrak{P}'^{(2)}\right)
\right)\right\rrbracket_{s}
\tag{1}
\\&=
\left\llbracket 
\mathrm{sc}^{0\mathbf{T}_{\Sigma^{\boldsymbol{\mathcal{A}}^{(2)}}}(X)}_{s}\left(
\mathrm{CH}^{(2)}_{s}\left(
\mathfrak{P}'^{(2)}
\right)
\right)
\right\rrbracket_{s}
\tag{2}
\\&=
\left\llbracket 
\mathrm{sc}^{0\mathbf{T}_{\Sigma^{\boldsymbol{\mathcal{A}}^{(2)}}}(X)}_{s}\left(
\mathrm{CH}^{(2)}_{s}\left(
\mathfrak{Q}'^{(2)}
\right)
\right)
\right\rrbracket_{s}
\tag{3}
\\&=
\left\llbracket \mathrm{CH}^{(2)}_{s}\left(
\mathrm{sc}^{0\mathbf{Pth}_{\boldsymbol{\mathcal{A}}^{(2)}}}_{s}\left(
\mathfrak{Q}'^{(2)}\right)
\right)\right\rrbracket_{s}
\tag{4}
\\&=
\left\llbracket \mathrm{CH}^{(2)}_{s}\left(
\mathfrak{Q}^{(2)}
\right)\right\rrbracket_{s}.
\tag{5}
\end{align*}

In the just stated chain of equalities, the first equality unravels the description of the second-order path $\mathfrak{P}^{(2)}$; the second equality follows from the fact that
\[
\left(
\mathrm{CH}^{(2)}_{s}\left(
\mathrm{sc}^{0\mathbf{Pth}_{\boldsymbol{\mathcal{A}}^{(2)}}}_{s}\left(
\mathfrak{P}'^{(2)}\right)
\right),
\mathrm{sc}^{0\mathbf{T}_{\Sigma^{\boldsymbol{\mathcal{A}}^{(2)}}}(X)}_{s}\left(
\mathrm{CH}^{(2)}_{s}\left(
\mathfrak{P}'^{(2)}
\right)
\right)
\right)\in \Theta^{(2)}_{s}
\]
and the fact that $\Theta^{(2)}\subseteq \Theta^{[2]}\subseteq \Theta^{\llbracket 2\rrbracket}$; the third equality follows from the fact that $\Theta^{\llbracket 2\rrbracket}$ is a $\Sigma^{\boldsymbol{\mathcal{A}}^{(2)}}$-congruence and by the fact that, by induction,  the following equality holds
\[
\left\llbracket
\mathrm{CH}^{(2)}_{s}\left(
\mathfrak{P}'^{(2)}
\right)
\right\rrbracket_{\Theta^{\llbracket 2\rrbracket}_{s}}
=
\left\llbracket
\mathrm{CH}^{(2)}_{s}\left(
\mathfrak{Q}'^{(2)}
\right)
\right\rrbracket_{\Theta^{\llbracket 2\rrbracket}_{s}};
\]
the fourth equality the second equality follows from the fact that
\[
\left(
\mathrm{CH}^{(2)}_{s}\left(
\mathrm{sc}^{0\mathbf{Pth}_{\boldsymbol{\mathcal{A}}^{(2)}}}_{s}\left(
\mathfrak{Q}'^{(2)}\right)
\right),
\mathrm{sc}^{0\mathbf{T}_{\Sigma^{\boldsymbol{\mathcal{A}}^{(2)}}}(X)}_{s}\left(
\mathrm{CH}^{(2)}_{s}\left(
\mathfrak{Q}'^{(2)}
\right)
\right)
\right)\in \Theta^{(2)}_{s}
\]
and the fact that $\Theta^{(2)}\subseteq \Theta^{[2]}\subseteq \Theta^{\llbracket 2\rrbracket}$; finally, the last equality recovers  the description of the second-order path $\mathfrak{Q}^{(2)}$.

This completes the case of the $0$-source.

\textsf{$\gamma$ is the $0$-target operation symbol.}

Let $(\mathfrak{P}'^{(2)}, \mathfrak{Q}'^{(2)})$ be the pair in $\mathrm{C}^{n}_{\boldsymbol{\mathcal{A}}^{(2)}}(\Upsilon^{(1)})_{s}$ for which
\[
\left(
\mathfrak{P}^{(2)},
\mathfrak{Q}^{(2)}
\right)
=
\left(
\mathrm{tg}^{0\mathbf{Pth}_{\boldsymbol{\mathcal{A}}^{(2)}}}_{s}\left(
\mathfrak{P}'^{(2)}\right),
\mathrm{tg}^{0\mathbf{Pth}_{\boldsymbol{\mathcal{A}}^{(2)}}}_{s}\left(
\mathfrak{Q}'^{(2)}\right)
\right).
\]

This case can be proven by a similar argument to that used for Case~(2.2).

This completes the case of the $0$-target.

\textsf{$\gamma$ is the $0$-composition operation symbol.}

Let $(\mathfrak{P}'^{(2)}, \mathfrak{Q}'^{(2)})$ and $(\mathfrak{P}''^{(2)}, \mathfrak{Q}''^{(2)})$ be the pairs in $\mathrm{C}^{n}_{\boldsymbol{\mathcal{A}}^{(2)}}(\Upsilon^{(1)})_{s}$ satisfying that
\begin{align*}
\mathrm{sc}^{(0,2)}_{s}\left(
\mathfrak{P}''^{(2)}
\right)
&=
\mathrm{tg}^{(0,2)}_{s}\left(
\mathfrak{P}'^{(2)}
\right);
&
\mathrm{sc}^{(0,2)}_{s}\left(
\mathfrak{Q}''^{(2)}
\right)
&=
\mathrm{tg}^{(0,2)}_{s}\left(
\mathfrak{Q}'^{(2)}
\right),
\end{align*}
for which
\[
\left(
\mathfrak{P}^{(2)},
\mathfrak{Q}^{(2)}
\right)
=
\left(
\mathfrak{P}''^{(2)}
\circ^{0\mathbf{Pth}_{\boldsymbol{\mathcal{A}}^{(2)}}}_{s}
\mathfrak{P}'^{(2)},
\mathfrak{Q}''^{(2)}
\circ^{0\mathbf{Pth}_{\boldsymbol{\mathcal{A}}^{(2)}}}_{s}
\mathfrak{Q}'^{(2)}
\right).
\]

Note that the following chain of equalities holds
\allowdisplaybreaks
\begin{align*}
\left\llbracket \mathrm{CH}^{(2)}_{s}\left(
\mathfrak{P}^{(2)}
\right)\right\rrbracket_{s}&=
\left\llbracket \mathrm{CH}^{(2)}_{s}\left(
\mathfrak{P}''^{(2)}
\circ^{0\mathbf{Pth}_{\boldsymbol{\mathcal{A}}^{(2)}}}_{s}
\mathfrak{P}'^{(2)}
\right)
\right\rrbracket_{s}
\tag{1}
\\&=
\left\llbracket 
\mathrm{CH}^{(2)}_{s}\left(
\mathfrak{P}''^{(2)}
\right)
\circ^{0\mathbf{T}_{\Sigma^{\boldsymbol{\mathcal{A}}^{(2)}}}(X)}_{s}
\mathrm{CH}^{(2)}_{s}\left(
\mathfrak{P}'^{(2)}
\right)
\right\rrbracket_{s}
\tag{2}
\\&=
\left\llbracket 
\mathrm{CH}^{(2)}_{s}\left(
\mathfrak{Q}''^{(2)}
\right)
\circ^{0\mathbf{T}_{\Sigma^{\boldsymbol{\mathcal{A}}^{(2)}}}(X)}_{s}
\mathrm{CH}^{(2)}_{s}\left(
\mathfrak{Q}'^{(2)}
\right)
\right\rrbracket_{s}
\tag{3}
\\&=
\left\llbracket \mathrm{CH}^{(2)}_{s}\left(
\mathfrak{P}''^{(2)}
\circ^{0\mathbf{Pth}_{\boldsymbol{\mathcal{A}}^{(2)}}}_{s}
\mathfrak{P}'^{(2)}
\right)\right\rrbracket_{s}
\tag{4}
\\&=
\left\llbracket \mathrm{CH}^{(2)}_{s}\left(
\mathfrak{Q}^{(2)}
\right)\right\rrbracket_{s}.
\tag{5}
\end{align*}

In the just stated chain of equalities, the first equality unravels the description of the second-order path $\mathfrak{P}^{(2)}$; the second equality follows from the fact that
\[
\left(
\mathrm{CH}^{(2)}_{s}\left(
\mathfrak{P}''^{(2)}
\circ^{0\mathbf{Pth}_{\boldsymbol{\mathcal{A}}^{(2)}}}_{s}
\mathfrak{P}'^{(2)}
\right),
\mathrm{CH}^{(2)}_{s}\left(
\mathfrak{P}''^{(2)}
\right)
\circ^{0\mathbf{T}_{\Sigma^{\boldsymbol{\mathcal{A}}^{(2)}}}(X)}_{s}
\mathrm{CH}^{(2)}_{s}\left(
\mathfrak{P}'^{(2)}
\right)
\right)\in \Theta^{(2)}_{s}
\]
and the fact that $\Theta^{(2)}\subseteq \Theta^{[2]}\subseteq \Theta^{\llbracket 2\rrbracket}$; the third equality follows from the fact that $\Theta^{\llbracket 2\rrbracket}$ is a $\Sigma^{\boldsymbol{\mathcal{A}}^{(2)}}$-congruence and by the fact that, by induction,  the following equalities holds
\allowdisplaybreaks
\begin{align*}
\left\llbracket
\mathrm{CH}^{(2)}_{s}\left(
\mathfrak{P}''^{(2)}
\right)
\right\rrbracket_{\Theta^{\llbracket 2\rrbracket}_{s}}
&=
\left\llbracket
\mathrm{CH}^{(2)}_{s}\left(
\mathfrak{Q}''^{(2)}
\right)
\right\rrbracket_{\Theta^{\llbracket 2\rrbracket}_{s}};
\\
\left\llbracket
\mathrm{CH}^{(2)}_{s}\left(
\mathfrak{P}'^{(2)}
\right)
\right\rrbracket_{\Theta^{\llbracket 2\rrbracket}_{s}}
&=
\left\llbracket
\mathrm{CH}^{(2)}_{s}\left(
\mathfrak{Q}'^{(2)}
\right)
\right\rrbracket_{\Theta^{\llbracket 2\rrbracket}_{s}};
\end{align*}
the fourth equality the second equality follows from the fact that
\[
\left(
\mathrm{CH}^{(2)}_{s}\left(
\mathfrak{Q}''^{(2)}
\circ^{0\mathbf{Pth}_{\boldsymbol{\mathcal{A}}^{(2)}}}_{s}
\mathfrak{Q}'^{(2)}
\right),
\mathrm{CH}^{(2)}_{s}\left(
\mathfrak{Q}''^{(2)}
\right)
\circ^{0\mathbf{T}_{\Sigma^{\boldsymbol{\mathcal{A}}^{(2)}}}(X)}_{s}
\mathrm{CH}^{(2)}_{s}\left(
\mathfrak{Q}'^{(2)}
\right)
\right)\in \Theta^{(2)}_{s}
\]
and the fact that $\Theta^{(2)}\subseteq \Theta^{[2]}\subseteq \Theta^{\llbracket 2\rrbracket}$; finally, the last equality recovers  the description of the second-order path $\mathfrak{Q}^{(2)}$.

This completes the case of the $0$-composition.

\textsf{$\gamma$ is the $1$-source operation symbol.}

Let $(\mathfrak{P}'^{(2)}, \mathfrak{Q}'^{(2)})$ be the pair in $\mathrm{C}^{n}_{\boldsymbol{\mathcal{A}}^{(2)}}(\Upsilon^{(1)})_{s}$ for which
\[
\left(
\mathfrak{P}^{(2)},
\mathfrak{Q}^{(2)}
\right)
=
\left(
\mathrm{sc}^{1\mathbf{Pth}_{\boldsymbol{\mathcal{A}}^{(2)}}}_{s}\left(
\mathfrak{P}'^{(2)}\right),
\mathrm{sc}^{1\mathbf{Pth}_{\boldsymbol{\mathcal{A}}^{(2)}}}_{s}\left(
\mathfrak{Q}'^{(2)}\right)
\right).
\]

Note that the following chain of equalities holds
\allowdisplaybreaks
\begin{align*}
\left\llbracket \mathrm{CH}^{(2)}_{s}\left(
\mathfrak{P}^{(2)}
\right)\right\rrbracket_{s}&=
\left\llbracket \mathrm{CH}^{(2)}_{s}\left(
\mathrm{sc}^{1\mathbf{Pth}_{\boldsymbol{\mathcal{A}}^{(2)}}}_{s}\left(
\mathfrak{P}'^{(2)}\right)
\right)\right\rrbracket_{s}
\tag{1}
\\&=
\left\llbracket 
\mathrm{sc}^{1\mathbf{T}_{\Sigma^{\boldsymbol{\mathcal{A}}^{(2)}}}(X)}_{s}\left(
\mathrm{CH}^{(2)}_{s}\left(
\mathfrak{P}'^{(2)}
\right)
\right)
\right\rrbracket_{s}
\tag{2}
\\&=
\left\llbracket 
\mathrm{sc}^{1\mathbf{T}_{\Sigma^{\boldsymbol{\mathcal{A}}^{(2)}}}(X)}_{s}\left(
\mathrm{CH}^{(2)}_{s}\left(
\mathfrak{Q}'^{(2)}
\right)
\right)
\right\rrbracket_{s}
\tag{3}
\\&=
\left\llbracket \mathrm{CH}^{(2)}_{s}\left(
\mathrm{sc}^{1\mathbf{Pth}_{\boldsymbol{\mathcal{A}}^{(2)}}}_{s}\left(
\mathfrak{Q}'^{(2)}\right)
\right)\right\rrbracket_{s}
\tag{4}
\\&=
\left\llbracket \mathrm{CH}^{(2)}_{s}\left(
\mathfrak{Q}^{(2)}
\right)\right\rrbracket_{s}.
\tag{5}
\end{align*}

In the just stated chain of equalities, the first equality unravels the description of the second-order path $\mathfrak{P}^{(2)}$; the second equality follows from the fact that
\[
\left(
\mathrm{CH}^{(2)}_{s}\left(
\mathrm{sc}^{1\mathbf{Pth}_{\boldsymbol{\mathcal{A}}^{(2)}}}_{s}\left(
\mathfrak{P}'^{(2)}\right)
\right),
\mathrm{sc}^{1\mathbf{T}_{\Sigma^{\boldsymbol{\mathcal{A}}^{(2)}}}(X)}_{s}\left(
\mathrm{CH}^{(2)}_{s}\left(
\mathfrak{P}'^{(2)}
\right)
\right)
\right)\in \Theta^{(2)}_{s}
\]
and the fact that $\Theta^{(2)}\subseteq \Theta^{[2]}\subseteq \Theta^{\llbracket 2\rrbracket}$; the third equality follows from the fact that $\Theta^{\llbracket 2\rrbracket}$ is a $\Sigma^{\boldsymbol{\mathcal{A}}^{(2)}}$-congruence and by the fact that, by induction,  the following equality holds
\[
\left\llbracket
\mathrm{CH}^{(2)}_{s}\left(
\mathfrak{P}'^{(2)}
\right)
\right\rrbracket_{\Theta^{\llbracket 2\rrbracket}_{s}}
=
\left\llbracket
\mathrm{CH}^{(2)}_{s}\left(
\mathfrak{Q}'^{(2)}
\right)
\right\rrbracket_{\Theta^{\llbracket 2\rrbracket}_{s}};
\]
the fourth equality the second equality follows from the fact that
\[
\left(
\mathrm{CH}^{(2)}_{s}\left(
\mathrm{sc}^{1\mathbf{Pth}_{\boldsymbol{\mathcal{A}}^{(2)}}}_{s}\left(
\mathfrak{Q}'^{(2)}\right)
\right),
\mathrm{sc}^{1\mathbf{T}_{\Sigma^{\boldsymbol{\mathcal{A}}^{(2)}}}(X)}_{s}\left(
\mathrm{CH}^{(2)}_{s}\left(
\mathfrak{Q}'^{(2)}
\right)
\right)
\right)\in \Theta^{(2)}_{s}
\]
and the fact that $\Theta^{(2)}\subseteq \Theta^{[2]}\subseteq \Theta^{\llbracket 2\rrbracket}$; finally, the last equality recovers  the description of the second-order path $\mathfrak{Q}^{(2)}$.

This completes the case of the $1$-source.

\textsf{$\gamma$ is the $1$-target operation symbol.}

Let $(\mathfrak{P}'^{(2)}, \mathfrak{Q}'^{(2)})$ be the pair in $\mathrm{C}^{n}_{\boldsymbol{\mathcal{A}}^{(2)}}(\Upsilon^{(1)})_{s}$ for which
\[
\left(
\mathfrak{P}^{(2)},
\mathfrak{Q}^{(2)}
\right)
=
\left(
\mathrm{tg}^{1\mathbf{Pth}_{\boldsymbol{\mathcal{A}}^{(2)}}}_{s}\left(
\mathfrak{P}'^{(2)}\right),
\mathrm{tg}^{1\mathbf{Pth}_{\boldsymbol{\mathcal{A}}^{(2)}}}_{s}\left(
\mathfrak{Q}'^{(2)}\right)
\right).
\]

This case can be proven by a similar argument to that used for Case~(2.5).

This completes the case of the $1$-target.

\textsf{$\gamma$ is the $1$-composition operation symbol.}

Let $(\mathfrak{P}'^{(2)}, \mathfrak{Q}'^{(2)})$ and $(\mathfrak{P}''^{(2)}, \mathfrak{Q}''^{(2)})$ be the pairs in $\mathrm{C}^{n}_{\boldsymbol{\mathcal{A}}^{(2)}}(\Upsilon^{(1)})_{s}$ satisfying that
\begin{align*}
\mathrm{sc}^{([1],2)}_{s}\left(
\mathfrak{P}''^{(2)}
\right)
&=
\mathrm{tg}^{([1],2)}_{s}\left(
\mathfrak{P}'^{(2)}
\right);
&
\mathrm{sc}^{([1],2)}_{s}\left(
\mathfrak{Q}''^{(2)}
\right)
&=
\mathrm{tg}^{([1],2)}_{s}\left(
\mathfrak{Q}'^{(2)}
\right),
\end{align*}
for which
\[
\left(
\mathfrak{P}^{(2)},
\mathfrak{Q}^{(2)}
\right)
=
\left(
\mathfrak{P}''^{(2)}
\circ^{1\mathbf{Pth}_{\boldsymbol{\mathcal{A}}^{(2)}}}_{s}
\mathfrak{P}'^{(2)},
\mathfrak{Q}''^{(2)}
\circ^{1\mathbf{Pth}_{\boldsymbol{\mathcal{A}}^{(2)}}}_{s}
\mathfrak{Q}'^{(2)}
\right).
\]

Note that the following chain of equalities holds
\allowdisplaybreaks
\begin{align*}
\left\llbracket \mathrm{CH}^{(2)}_{s}\left(
\mathfrak{P}^{(2)}
\right)\right\rrbracket_{s}&=
\left\llbracket \mathrm{CH}^{(2)}_{s}\left(
\mathfrak{P}''^{(2)}
\circ^{1\mathbf{Pth}_{\boldsymbol{\mathcal{A}}^{(2)}}}_{s}
\mathfrak{P}'^{(2)}
\right)
\right\rrbracket_{s}
\tag{1}
\\&=
\left\llbracket 
\mathrm{CH}^{(2)}_{s}\left(
\mathfrak{P}''^{(2)}
\right)
\circ^{1\mathbf{T}_{\Sigma^{\boldsymbol{\mathcal{A}}^{(2)}}}(X)}_{s}
\mathrm{CH}^{(2)}_{s}\left(
\mathfrak{P}'^{(2)}
\right)
\right\rrbracket_{s}
\tag{2}
\\&=
\left\llbracket 
\mathrm{CH}^{(2)}_{s}\left(
\mathfrak{Q}''^{(2)}
\right)
\circ^{1\mathbf{T}_{\Sigma^{\boldsymbol{\mathcal{A}}^{(2)}}}(X)}_{s}
\mathrm{CH}^{(2)}_{s}\left(
\mathfrak{Q}'^{(2)}
\right)
\right\rrbracket_{s}
\tag{3}
\\&=
\left\llbracket \mathrm{CH}^{(2)}_{s}\left(
\mathfrak{P}''^{(2)}
\circ^{1\mathbf{Pth}_{\boldsymbol{\mathcal{A}}^{(2)}}}_{s}
\mathfrak{P}'^{(2)}
\right)\right\rrbracket_{s}
\tag{4}
\\&=
\left\llbracket \mathrm{CH}^{(2)}_{s}\left(
\mathfrak{Q}^{(2)}
\right)\right\rrbracket_{s}.
\tag{5}
\end{align*}

In the just stated chain of equalities, the first equality unravels the description of the second-order path $\mathfrak{P}^{(2)}$; the second equality follows from the fact that
\[
\left(
\mathrm{CH}^{(2)}_{s}\left(
\mathfrak{P}''^{(2)}
\circ^{1\mathbf{Pth}_{\boldsymbol{\mathcal{A}}^{(2)}}}_{s}
\mathfrak{P}'^{(2)}
\right),
\mathrm{CH}^{(2)}_{s}\left(
\mathfrak{P}''^{(2)}
\right)
\circ^{1\mathbf{T}_{\Sigma^{\boldsymbol{\mathcal{A}}^{(2)}}}(X)}_{s}
\mathrm{CH}^{(2)}_{s}\left(
\mathfrak{P}'^{(2)}
\right)
\right)\in \Theta^{(2)}_{s}
\]
and the fact that $\Theta^{(2)}\subseteq \Theta^{[2]}\subseteq \Theta^{\llbracket 2\rrbracket}$; the third equality follows from the fact that $\Theta^{\llbracket 2\rrbracket}$ is a $\Sigma^{\boldsymbol{\mathcal{A}}^{(2)}}$-congruence and by the fact that, by induction,  the following equalities holds
\allowdisplaybreaks
\begin{align*}
\left\llbracket
\mathrm{CH}^{(2)}_{s}\left(
\mathfrak{P}''^{(2)}
\right)
\right\rrbracket_{\Theta^{\llbracket 2\rrbracket}_{s}}
&=
\left\llbracket
\mathrm{CH}^{(2)}_{s}\left(
\mathfrak{Q}''^{(2)}
\right)
\right\rrbracket_{\Theta^{\llbracket 2\rrbracket}_{s}};
\\
\left\llbracket
\mathrm{CH}^{(2)}_{s}\left(
\mathfrak{P}'^{(2)}
\right)
\right\rrbracket_{\Theta^{\llbracket 2\rrbracket}_{s}}
&=
\left\llbracket
\mathrm{CH}^{(2)}_{s}\left(
\mathfrak{Q}'^{(2)}
\right)
\right\rrbracket_{\Theta^{\llbracket 2\rrbracket}_{s}};
\end{align*}
the fourth equality the second equality follows from the fact that
\[
\left(
\mathrm{CH}^{(2)}_{s}\left(
\mathfrak{Q}''^{(2)}
\circ^{1\mathbf{Pth}_{\boldsymbol{\mathcal{A}}^{(2)}}}_{s}
\mathfrak{Q}'^{(2)}
\right),
\mathrm{CH}^{(2)}_{s}\left(
\mathfrak{Q}''^{(2)}
\right)
\circ^{1\mathbf{T}_{\Sigma^{\boldsymbol{\mathcal{A}}^{(2)}}}(X)}_{s}
\mathrm{CH}^{(2)}_{s}\left(
\mathfrak{Q}'^{(2)}
\right)
\right)\in \Theta^{(2)}_{s}
\]
and the fact that $\Theta^{(2)}\subseteq \Theta^{[2]}\subseteq \Theta^{\llbracket 2\rrbracket}$; finally, the last equality recovers  the description of the second-order path $\mathfrak{Q}^{(2)}$.

This completes the case of the $1$-composition.

This completes Case~(2).

This completes the proof.
\end{proof}

In the following proposition we will prove that, for every sort $s\in S$, if two second-order paths in $\mathrm{Pth}_{\boldsymbol{\mathcal{A}}^{(2)},s}$ are $\mathrm{Ker}(\mathrm{CH}^{(2)})_{s}\vee\Upsilon^{[1]}_{s}$-related, then their respective values under the second-order Curry-Howard mapping are $\Theta^{\llbracket 2\rrbracket}_{s}$-related.

\begin{restatable}{proposition}{PDVDCH}
\label{PDVDCH} Let $s$ be a sort in $S$ and let $\mathfrak{P}^{(2)}$ and $\mathfrak{Q}^{(2)}$ be second-order paths in $\mathrm{Pth}_{\boldsymbol{\mathcal{A}}^{(2)},s}$ satisfying that $\llbracket\mathfrak{P}^{(2)}\rrbracket_{s}=\llbracket\mathfrak{Q}^{(2)}\rrbracket_{s}$, then
\[
\Bigl\llbracket
\mathrm{CH}^{(2)}_{s}\left(
\mathfrak{P}^{(2)}
\right)
\Bigr\rrbracket_{\Theta^{\llbracket 2\rrbracket}_{s}}
=
\Bigl\llbracket
\mathrm{CH}^{(2)}_{s}\left(
\mathfrak{Q}^{(2)}
\right)
\Bigr\rrbracket_{\Theta^{\llbracket 2\rrbracket}_{s}}.
\]
\end{restatable}
\begin{proof}
Let us recall that, for a sort $s\in S$, two second-order paths $\mathfrak{P}^{(2)}$ and $\mathfrak{Q}^{(2)}$ in $\mathrm{Pth}_{\boldsymbol{\mathcal{A}}^{(2)},s}$ are related under $\mathrm{Ker}(\mathrm{CH}^{(2)}_{s})\vee \Upsilon^{[1]}_{s}$ if, and only if, there exists a natural number $m\in \mathbb{N}-\{0\}$ and a sequence of second-order paths $(\mathfrak{R}^{(2)}_{k})_{k\in m+1}$ in $\mathrm{Pth}_{\boldsymbol{\mathcal{A}}^{(2)},s}^{m+1}$ satisfying that 
\begin{enumerate}
\item $\mathfrak{R}^{(2)}_{0}=\mathfrak{P}^{(2)}$;
\item $\mathfrak{R}^{(2)}_{m}=\mathfrak{Q}^{(2)}$;
\item for every $k\in m$, 
$(\mathfrak{R}^{(2)}_{k},\mathfrak{R}^{(2)}_{k+1})\in 
\mathrm{Ker}(\mathrm{CH}^{(2)})_{s}\cup \Upsilon^{[1]}_{s}$.
\end{enumerate}

We prove the statement by induction on $m$.

\textsf{Base step of the induction.}

For the case $m=1$, we have that either the pair $(\mathfrak{P}^{(2)},\mathfrak{Q}^{(2)})$ is in $\mathrm{Ker}(\mathrm{CH}^{(2)})_{s}$ or in $\Upsilon^{[1]}_{s}$.  

In case $(\mathfrak{P}^{(2)},\mathfrak{Q}^{(2)})$ is in $\mathrm{Ker}(\mathrm{CH}^{(2)})_{s}$, then we have that $\mathrm{CH}^{(2)}_{s}(\mathfrak{P}^{(2)})=\mathrm{CH}^{(2)}_{s}(\mathfrak{Q}^{(2)})$. This case follows by reflexivity of $\Theta^{\llbracket 2 \rrbracket}_{s}$.

In case $(\mathfrak{P}^{(2)},\mathfrak{Q}^{(2)})$ is in $\Upsilon^{[1]}_{s}$, the statement holds by Proposition~\ref{PDUpsDCH}.

\textsf{Inductive step of the induction.}

Assume that the statement holds for sequences of length $m$, i.e., if $(\mathfrak{R}^{(2)}_{k})_{k\in m+1}$ is a sequence of length $m$ in $\mathrm{Pth}_{\boldsymbol{\mathcal{A}}^{(2)},s}^{m+1}$ satisfying that 
\begin{enumerate}
\item $\mathfrak{R}^{(2)}_{0}=\mathfrak{P}^{(2)}$;
\item $\mathfrak{R}^{(2)}_{m}=\mathfrak{Q}^{(2)}$;
\item for every $k\in m$, 
$(\mathfrak{R}^{(2)}_{k},\mathfrak{R}^{(2)}_{k+1})\in 
\mathrm{Ker}(\mathrm{CH}^{(2)})_{s}\cup \Upsilon^{[1]}_{s}$,
\end{enumerate}
then it holds that
\[
\Bigl\llbracket
\mathrm{CH}^{(2)}_{s}\left(
\mathfrak{P}^{(2)}
\right)
\Bigr\rrbracket_{\Theta^{\llbracket 2\rrbracket}_{s}}
=
\Bigl\llbracket
\mathrm{CH}^{(2)}_{s}\left(
\mathfrak{Q}^{(2)}
\right)
\Bigr\rrbracket_{\Theta^{\llbracket 2\rrbracket}_{s}}.
\]

Now we prove it for a sequence of length $m+1$. Let $(\mathfrak{R}^{(2)}_{k})_{k\in m+2}$ is a sequence of length $m+1$ in $\mathrm{Pth}_{\boldsymbol{\mathcal{A}}^{(2)},s}^{m+2}$ satisfying that 
\begin{enumerate}
\item $\mathfrak{R}^{(2)}_{0}=\mathfrak{P}^{(2)}$;
\item $\mathfrak{R}^{(2)}_{m+1}=\mathfrak{Q}^{(2)}$;
\item for every $k\in m+1$, 
$(\mathfrak{R}^{(2)}_{k},\mathfrak{R}^{(2)}_{k+1})\in 
\mathrm{Ker}(\mathrm{CH}^{(2)})_{s}\cup \Upsilon^{[1]}_{s}$.
\end{enumerate}

Consider the initial subsequence $(\mathfrak{R}^{(2)}_{k})_{k\in m+1}$. It is a sequence of length $m$ in $\mathrm{Pth}_{\boldsymbol{\mathcal{A}}^{(2)},s}^{m+1}$ satisfying that 
\begin{enumerate}
\item $\mathfrak{R}^{(2)}_{0}=\mathfrak{P}^{(2)}$;
\item $\mathfrak{R}^{(2)}_{m}=\mathfrak{R}^{(2)}_{m}$;
\item for every $k\in m$, 
$(\mathfrak{R}^{(2)}_{k},\mathfrak{R}^{(2)}_{k+1})\in 
\mathrm{Ker}(\mathrm{CH}^{(2)})_{s}\cup \Theta^{[1]}_{s}$.
\end{enumerate}
By the inductive hypothesis it holds that 
\[
\Bigl\llbracket
\mathrm{CH}^{(2)}_{s}\left(
\mathfrak{P}^{(2)}
\right)
\Bigr\rrbracket_{\Theta^{\llbracket 2\rrbracket}_{s}}
=
\Bigl\llbracket
\mathrm{CH}^{(2)}_{s}\left(
\mathfrak{R}^{(2)}_{m}
\right)
\Bigr\rrbracket_{\Theta^{\llbracket 2\rrbracket}_{s}}.
\]

Now consider the final subsequence $(\mathfrak{R}^{(2)}_{m},\mathfrak{Q}^{(2)})$. It is a sequence of length $1$ in $\mathrm{Pth}_{\boldsymbol{\mathcal{A}}^{(2)},s}^{2}$ from $\mathfrak{R}^{(2)}_{m}$ to $\mathfrak{Q}^{(2)}$ satisfying that 
$(\mathfrak{R}^{(2)}_{m},\mathfrak{Q}^{(2)})\in 
\mathrm{Ker}(\mathrm{CH}^{(2)})_{s}\cup \Theta^{[1]}_{s}$. By the base case, we have that 
\[
\Bigl\llbracket
\mathrm{CH}^{(2)}_{s}\left(
\mathfrak{R}^{(2)}_{m}
\right)
\Bigr\rrbracket_{\Theta^{\llbracket 2\rrbracket}_{s}}
=
\Bigl\llbracket
\mathrm{CH}^{(2)}_{s}\left(
\mathfrak{Q}^{(2)}
\right)
\Bigr\rrbracket_{\Theta^{\llbracket 2\rrbracket}_{s}}.
\]

By transitivity of $\Theta^{\llbracket 2\rrbracket}_{s}$ we conclude that 
\[
\Bigl\llbracket
\mathrm{CH}^{(2)}_{s}\left(
\mathfrak{P}^{(2)}
\right)
\Bigr\rrbracket_{\Theta^{\llbracket 2\rrbracket}_{s}}
=
\Bigl\llbracket
\mathrm{CH}^{(2)}_{s}\left(
\mathfrak{Q}^{(2)}
\right)
\Bigr\rrbracket_{\Theta^{\llbracket 2\rrbracket}_{s}}.
\]

This completes the proof.
\end{proof}

The following lemma states that if a term of sort $s$ is such that its image under $\mathrm{ip}^{(2,X)@}$ is a second-order path, then it is related with respect to the $\Sigma^{\boldsymbol{\mathcal{A}}^{(2)}}$-congruence $\Theta^{\llbracket 2\rrbracket }$ on $\mathbf{T}_{\Sigma^{\boldsymbol{\mathcal{A}}^{(2)}}}(X)$, with a term of the same sort in $\mathrm{CH}^{(2)}[\mathrm{Pth}_{\boldsymbol{\mathcal{A}}^{(2)}}]$. Actually, we prove that such a term is related with its image under the action of $\mathrm{CH}^{(2)}\circ \mathrm{ip}^{(2,X)@}$.

\begin{restatable}{lemma}{LDVWCong}
\label{LDVWCong} Let $s$ be a sort in $S$ and $P$ a term in $\mathrm{T}_{\Sigma^{\boldsymbol{\mathcal{A}}}}(X)_{s}$. If $\mathrm{ip}^{(2,X)@}_{s}(P)$ is a second-order path in $\mathrm{Pth}_{\boldsymbol{\mathcal{A}}^{(2)},s}$ then $(P,\mathrm{CH}^{(2)}_{s}(\mathrm{ip}^{(2,X)@}_{s}(P)))\in\Theta^{\llbracket 2\rrbracket}_{s}$.
\end{restatable}
\begin{proof}
Since $\mathrm{ip}^{(2,X)@}_{s}(P)$ is a second-order path in $\mathrm{Pth}_{\boldsymbol{\mathcal{A}}^{(2)},s}$, we conclude according to Lemma~\ref{LDWCong} that
\[
\left(
P,\,
\mathrm{CH}^{(2)}_{s}\left(
\mathrm{ip}^{(2,X)@}_{s}\left(
P
\right)\right)\right)\in\Theta^{[2]}_{s}.
\]
From the above statement, we conclude that the pair above is also $\Theta^{\llbracket 2\rrbracket}_{s}$-related.

This concludes the proof.
\end{proof}

We continue with the description of the behaviour of the $\Sigma^{\boldsymbol{\mathcal{A}}^{(2)}}$-congruence $\Theta^{\llbracket 2 \rrbracket}$. In the following lemma we state that, for every sort $s\in S$, if two terms are $\Theta^{\llbracket 2 \rrbracket}_{s}$-related and one of them, when mapped under $\mathrm{ip}^{(2,X)@}$, retrieves a second-order path, then the other term has a similar behaviour. Moreover, if the just described situation happens, then these two paths will be in the same $\llbracket \cdot \rrbracket$-class.

\begin{restatable}{lemma}{LDVCong}
\label{LDVCong} Let $s$ be a sort in $S$ and $P,Q$ terms in $\mathrm{T}_{\Sigma^{\boldsymbol{\mathcal{A}}^{(2)}}}(X)_{s}$ such that $(P,Q)\in\Theta^{\llbracket 2 \rrbracket}_{s}$, then
\begin{enumerate}
\item[(i)] $\mathrm{ip}^{(2,X)@}_{s}(P)\in\mathrm{Pth}_{\boldsymbol{\mathcal{A}}^{(2)},s}$ if, and only if, $\mathrm{ip}^{(2,X)@}_{s}(Q)\in\mathrm{Pth}_{\boldsymbol{\mathcal{A}}^{(2)},s}$;
\item[(ii)] If $\mathrm{ip}^{(2,X)@}_{s}(P)$ or $\mathrm{ip}^{(2,X)@}_{s}(Q)$ is a second-order path in $\mathrm{Pth}_{\boldsymbol{\mathcal{A}}^{(2)},s}$ then
\[
\Bigl\llbracket \mathrm{ip}^{(2,X)@}_{s}(P) \Bigr\rrbracket_{s}=
\Bigl\llbracket \mathrm{ip}^{(2,X)@}_{s}(Q) \Bigr\rrbracket_{s}.
\]
\end{enumerate}
\end{restatable}
\begin{proof}
Let us recall that, for a sort $s\in S$, two terms $P$ and $Q$  in $\mathrm{T}_{\Sigma^{\boldsymbol{\mathcal{A}}^{(2)}}}(X)_{s}$ are $\Theta^{\llbracket 2\rrbracket}_{s} $ if, and only if, there exists a natural number $m\in\mathbb{N}-\{0\}$ and a sequence of terms $(R_{k})_{k\in m+1}$ in $\mathrm{T}_{\Sigma^{\boldsymbol{\mathcal{A}}^{(2)}}}(X)_{s}^{m+1}$ satisfying that  
\begin{enumerate}
\item $R_{0}=P$;
\item $R_{m}=Q$;
\item for every $k\in m$, 
$(R_{k},R_{k+1})\in 
 \Theta^{[2]}_{s}\cup \Psi^{[1]}_{s}$.
\end{enumerate}

We prove the statement by induction on $m$.

\textsf{Base step of the induction.}

For the case $m=1$, we have that either the pair $(P,Q)$ is in $\Theta^{[2]}_{s}$ or in $\Psi^{[1]}_{s}$. The statement follows by Lemmas~\ref{LDThetaCong} and~\ref{LDPsiCong}, respectively.

\textsf{Inductive step of the induction.}

Assume that the statement holds for sequences of length $m$, i.e., if $(R_{k})_{k\in m+1}$ is a sequence of length $m$ in $\mathrm{T}_{\Sigma^{\boldsymbol{\mathcal{A}}^{(2)}}}(X)_{s}^{m+1}$ satisfying that 
\begin{enumerate}
\item $R_{0}=P$;
\item $R_{m}=Q$;
\item for every $k\in m$, 
$(R_{k},R_{k+1})\in 
 \Theta^{[2]}_{s}\cup \Psi^{[1]}_{s}$,
\end{enumerate}
then it holds that 
\begin{enumerate}
\item[(i)] $\mathrm{ip}^{(2,X)@}_{s}(P)\in\mathrm{Pth}_{\boldsymbol{\mathcal{A}}^{(2)},s}$ if, and only if, $\mathrm{ip}^{(2,X)@}_{s}(Q)\in\mathrm{Pth}_{\boldsymbol{\mathcal{A}}^{(2)},s}$;
\item[(ii)] If $\mathrm{ip}^{(2,X)@}_{s}(P)$ or $\mathrm{ip}^{(2,X)@}_{s}(Q)$ is a second-order path in $\mathrm{Pth}_{\boldsymbol{\mathcal{A}}^{(2)},s}$ then
\[
\Bigl\llbracket \mathrm{ip}^{(2,X)@}_{s}(P) \Bigr\rrbracket_{s}=
\Bigl\llbracket \mathrm{ip}^{(2,X)@}_{s}(Q) \Bigr\rrbracket_{s}.
\]
\end{enumerate}

Now we prove it for a sequence of length $m+1$. Let $(R_{k})_{k\in m+2}$ be a sequence of length $m+1$ in $\mathrm{T}_{\Sigma^{\boldsymbol{\mathcal{A}}^{(2)}}}(X)_{s}^{m+2}$ satisfying that 
\begin{enumerate}
\item $R_{0}=P$;
\item $R_{m+1}=Q$;
\item for every $k\in m+1$, 
$(R_{k},R_{k+1})\in 
 \Theta^{[2]}_{s}\cup \Psi^{[1]}_{s}$,
\end{enumerate}

Consider the initial subsequence $(R_{k})_{k\in m+1}$. It is a sequence of length $m$ in $\mathrm{T}_{\Sigma^{\boldsymbol{\mathcal{A}}^{(2)}}}(X)_{s}^{m+1}$ satisfying that 
\begin{enumerate}
\item $R_{0}=P$;
\item $R_{m}=R_{m}$;
\item for every $k\in m$, 
$(R_{k},R_{k+1})\in 
 \Theta^{[2]}_{s}\cup \Psi^{[1]}_{s}$.
\end{enumerate}

By the inductive hypothesis it holds that 
\begin{enumerate}
\item[(i)] $\mathrm{ip}^{(2,X)@}_{s}(P)\in\mathrm{Pth}_{\boldsymbol{\mathcal{A}}^{(2)},s}$ if, and only if, $\mathrm{ip}^{(2,X)@}_{s}(R_{m})\in\mathrm{Pth}_{\boldsymbol{\mathcal{A}}^{(2)},s}$;
\item[(ii)] If $\mathrm{ip}^{(2,X)@}_{s}(P)$ or $\mathrm{ip}^{(2,X)@}_{s}(R_{m})$ is a second-order path in $\mathrm{Pth}_{\boldsymbol{\mathcal{A}}^{(2)},s}$ then
\[
\Bigl\llbracket \mathrm{ip}^{(2,X)@}_{s}(P) \Bigr\rrbracket_{s}=
\Bigl\llbracket \mathrm{ip}^{(2,X)@}_{s}(R_{m}) \Bigr\rrbracket_{s}.
\]
\end{enumerate}

Now consider the final subsequence $(R_{m},Q)$. It is a sequence of length $1$ in $\mathrm{T}_{\Sigma^{\boldsymbol{\mathcal{A}}^{(2)}}}(X)_{s}^{2}$ from $R_{m}$ to $Q$ satisfying that  $(R_{m},Q)\in\Theta^{[2]}_{s}\cup \Psi^{[1]}_{s}$. By the base case, we have that 
\begin{enumerate}
\item[(i)] $\mathrm{ip}^{(2,X)@}_{s}(R_{m})\in\mathrm{Pth}_{\boldsymbol{\mathcal{A}}^{(2)},s}$ if, and only if, $\mathrm{ip}^{(2,X)@}_{s}(Q)\in\mathrm{Pth}_{\boldsymbol{\mathcal{A}}^{(2)},s}$;
\item[(ii)] If $\mathrm{ip}^{(2,X)@}_{s}(R_{m})$ or $\mathrm{ip}^{(2,X)@}_{s}(Q)$ is a second-order path in $\mathrm{Pth}_{\boldsymbol{\mathcal{A}}^{(2)},s}$ then
\[
\Bigl\llbracket \mathrm{ip}^{(2,X)@}_{s}(R_{m}) \Bigr\rrbracket_{s}=
\Bigl\llbracket \mathrm{ip}^{(2,X)@}_{s}(Q) \Bigr\rrbracket_{s}.
\]
\end{enumerate}

All in all, we conclude that 
\begin{enumerate}
\item[(i)] $\mathrm{ip}^{(2,X)@}_{s}(P)\in\mathrm{Pth}_{\boldsymbol{\mathcal{A}}^{(2)},s}$ if, and only if, $\mathrm{ip}^{(2,X)@}_{s}(Q)\in\mathrm{Pth}_{\boldsymbol{\mathcal{A}}^{(2)},s}$;
\item[(ii)] If $\mathrm{ip}^{(2,X)@}_{s}(P$ or $\mathrm{ip}^{(2,X)@}_{s}(Q)$ is a second-order path in $\mathrm{Pth}_{\boldsymbol{\mathcal{A}}^{(2)},s}$ then
\[
\Bigl\llbracket \mathrm{ip}^{(2,X)@}_{s}(P) \Bigr\rrbracket_{s}=
\Bigl\llbracket \mathrm{ip}^{(2,X)@}_{s}(Q) \Bigr\rrbracket_{s}.
\]
\end{enumerate}

This completes the proof.
\end{proof}

\begin{restatable}{corollary}{CDVCong}
\label{CDVCong} Let $s$ be a sort in $S$, $\mathfrak{P}^{(2)}$ a second-order path in $\mathrm{Pth}_{\boldsymbol{\mathcal{A}}^{(2)},s}$ and $P$ a term in $\mathrm{T}_{\Sigma^{\boldsymbol{\mathcal{A}}^{(2)}}}(X)_{s}$ such that $(P,\mathrm{CH}^{(2)}_{s}(\mathfrak{P}^{(2)}))\in\Theta^{\llbracket 2 \rrbracket}_{s}$. Then $\mathrm{ip}^{(2,X)@}_{s}(P)$ is a second-order path in $\llbracket \mathfrak{P}^{(2)}\rrbracket_{s}$.
\end{restatable}
\begin{proof}
Let us recall from Proposition~\ref{PDIpDCH} that $\mathrm{ip}^{(2,X)@}_{s}(\mathrm{CH}^{(2)}_{s}(\mathfrak{P}^{(2)}))$ is a second-order in $[\mathfrak{P}^{(2)}]_{s}$. Consequently, the $\llbracket \cdot \rrbracket$-class of the second-order paths $\mathfrak{P}^{(2)}$ and $\mathrm{ip}^{(2,X)@}_{s}(\mathrm{CH}^{(2)}_{s}(\mathfrak{P}^{(2)}))$ is the same.

According to Lemma~\ref{LDVCong} we have that $\mathrm{ip}^{(2,X)@}_{s}(P)$ is a second-order path in $\mathrm{Pth}_{\boldsymbol{\mathcal{A}}^{(2)},s}$. Moreover, we have that 
\[
\Bigl\llbracket \mathrm{ip}^{(2,X)@}_{s}(P) \Bigr\rrbracket_{s}=
\Bigl\llbracket \mathrm{ip}^{(2,X)@}_{s}\left(
\mathrm{CH}^{(2)}_{s}\left(
\mathfrak{P}^{(2)}
\right)
\right) \Bigr\rrbracket_{s}
=
\Bigl\llbracket \mathfrak{P}^{(2)}\Bigr\rrbracket_{s}
.
\]

This completes the proof.
\end{proof}

\begin{restatable}{corollary}{CDVCongDCH}
\label{CDVCongDCH} Let $s$ be a sort in $S$ and $\mathfrak{P}'^{(2)}$, $\mathfrak{P}^{(2)}$ two second-order paths in $\mathrm{Pth}_{\boldsymbol{\mathcal{A}}^{(2)},s}$. If $(\mathrm{CH}^{(2)}_{s}(\mathfrak{P}'^{(2)}),\mathrm{CH}^{(2)}_{s}(\mathfrak{P}^{(2)}))\in\Theta^{\llbracket 2 \rrbracket}_{s}$, then $\llbracket \mathfrak{P}'^{(2)}\rrbracket_{s}=\llbracket \mathfrak{P}^{(2)}\rrbracket_{s}$.
\end{restatable}
\begin{proof}
Following Corollary~\ref{CDVCong}, we have that $\mathrm{ip}^{(2,X)@}_{s}(\mathrm{CH}^{(2)}_{s}(\mathfrak{P}'^{(2)}))$ is a second-order path in $\llbracket \mathfrak{P}^{(2)}\rrbracket_{s}$.

However, according to Proposition~\ref{PDIpDCH}, $\mathrm{ip}^{(2,X)@}_{s}(\mathrm{CH}^{(2)}_{s}(\mathfrak{P}'^{(2)}))$ is a second-order in $[\mathfrak{P}'^{(2)}]_{s}$. Consequently, the $\llbracket \cdot \rrbracket$-class of the second-order paths  $\mathrm{ip}^{(2,X)@}_{s}(\mathrm{CH}^{(2)}_{s}(\mathfrak{P}'^{(2)}))$  and $\mathfrak{P}'^{(2)}$ is the same.

It follows that $\llbracket \mathfrak{P}'^{(2)}\rrbracket_{s}=\llbracket \mathfrak{P}^{(2)}\rrbracket_{s}$.
\end{proof}

\begin{remark}\label{RDVExplFalla} Taking into account Lemmas~\ref{LDVWCong} and~\ref{LDVCong} and Corollary~\ref{CDVCongDCH} we infer that, for every sort $s\in S$, if a term $P\in\mathrm{T}_{\Sigma^{\boldsymbol{\mathcal{A}}^{(2)}}}(X)_{s}$ belongs to $\llbracket \mathrm{CH}^{(2)}(\mathfrak{P}^{(2)})\rrbracket_{\Theta^{\llbracket 2 \rrbracket}_{s}}$, for some second-order path $\mathfrak{P}^{(2)}\in\mathrm{Pth}_{\boldsymbol{\mathcal{A}}^{(2)},s}$ then $P$ can be understood as an alternative term describing the same second-order path as $\mathfrak{P}^{(2)}$ up to equivalence with respect to $\mathrm{Ker}(\mathrm{CH}^{(2)})\vee \Upsilon^{[1]}$. Therefore, we can think of the class $\llbracket \mathrm{CH}^{(2)}_{s}(\mathfrak{P}^{(2)})\rrbracket_{\Theta^{\llbracket 2\rrbracket}_{s}}$ as a class containing all possible terms describing all second-order paths in $\llbracket\mathfrak{P}^{(2)}\rrbracket_{s}$. Given one such term, $P\in \llbracket\mathrm{CH}^{(2)}_{s}(\mathfrak{P}^{(2)})\rrbracket_{\Theta^{\llbracket 2\rrbracket}_{s}}$, then $\mathrm{ip}^{(2,X)@}_{s}(P)$ is a second-order path in $\llbracket\mathfrak{P}^{(2)}\rrbracket^{}_{s}$. 
\end{remark}

\chapter{
\texorpdfstring
{The $S$-sorted set $\mathrm{PT}_{\boldsymbol{\mathcal{A}}^{(2)}}$ of second-order path terms}
{Second-order path terms}
}\label{S2M}

In this chapter we introduce the many-sorted set of second-order path terms, denoted by  $\mathrm{PT}_{\boldsymbol{\mathcal{A}}^{(2)}}$, as the $\Theta^{\llbracket 2 \rrbracket}$-saturation of the subset $\mathrm{CH}^{(2)}[\mathrm{Pth}_{\boldsymbol{\mathcal{A}}^{(2)}}]$ of $\mathrm{T}_{\Sigma^{\boldsymbol{\mathcal{A}}^{(2)}}}(X)$. We show that a term $P$ in $\mathrm{T}_{\Sigma^{\boldsymbol{\mathcal{A}}^{(2)}}}(X)_{s}$  is a second-order path term if, and only if, $\mathrm{ip}^{(2,X)@}_{s}(P)$ is a second-order path in $\mathrm{Pth}_{\boldsymbol{\mathcal{A}}^{(2)},s}$. In fact this allows us to consider the corestriction of $\mathrm{ip}^{(2,X)@}$ to $\mathrm{Pth}_{\boldsymbol{\mathcal{A}}^{(2)}}$ and the correstriction of $\mathrm{CH}^{(2)}$ to $\mathrm{PT}_{\boldsymbol{\mathcal{A}}^{(2)}}$, that will keep the same notation. We prove that $\mathrm{PT}_{\boldsymbol{\mathcal{A}}^{(2)}}$ has an structure of many-sorted partial  $\Sigma^{\boldsymbol{\mathcal{A}}^{(2)}}$-algebra, denoted by $\mathbf{PT}_{\boldsymbol{\mathcal{A}}^{(2)}}$, inherited from $\mathbf{T}_{\Sigma^{\boldsymbol{\mathcal{A}}^{(2)}}}(X)/{\Theta^{\llbracket 2 \rrbracket}}$. We can prove that $\mathbf{PT}_{\boldsymbol{\mathcal{A}}^{(2)}}$ is a partial Dedekind-Peano $\Sigma^{\boldsymbol{\mathcal{A}}^{(2)}}$-algebra. This implies, in particular, that subterms of second-order path terms are itself second-order path terms. This implies that the many-sorted set  $\mathrm{PT}_{\boldsymbol{\mathcal{A}}^{(2)}}$ of second-order path terms also admits an Artinian partial order, denoted by $\leq_{\mathbf{PT}_{\boldsymbol{\mathcal{A}}^{(2)}}}(X)$, given by the standard subterm relation. We then introduce the quotient of second-order path terms by the restriction of the relation $\Theta^{\llbracket 2 \rrbracket}$, that we will denote by $\llbracket \mathrm{PT}_{\boldsymbol{\mathcal{A}}^{(2)}}\rrbracket$. The possible compositions of $\mathrm{ip}^{(2,X)@}$ and  $\mathrm{CH}^{(2)}$ with $\mathrm{pr}^{\Theta^{\llbracket 2 \rrbracket}}$, the projection to the quotient, will be denoted by $\mathrm{ip}^{(\llbracket 2 \rrbracket,X)@}$ and $\mathrm{CH}^{\llbracket 2 \rrbracket}$, respectively. It is shown that $\llbracket \mathrm{PT}_{\boldsymbol{\mathcal{A}}^{(2)}}\rrbracket$ has a structure of many-sorted partial  $\Sigma^{\boldsymbol{\mathcal{A}}^{(2)}}$-algebra, denoted $\llbracket \mathbf{PT}_{\boldsymbol{\mathcal{A}}^{(2)}}\rrbracket$. Moreover, $\llbracket \mathrm{PT}_{\boldsymbol{\mathcal{A}}^{(2)}}\rrbracket$ is an $S$-sorted $2$-category, more precisely, it has a structure of $2$-categorial $\Sigma$-algebra, denoted by  $\llbracket \mathsf{PT}_{\boldsymbol{\mathcal{A}}^{(2)}}\rrbracket$. Finally, the quotient $\llbracket \mathrm{PT}_{\boldsymbol{\mathcal{A}}^{(2)}}\rrbracket$ admits an Artinian partial preorder, denoted by  $\leq_{\llbracket \mathbf{PT}_{\boldsymbol{\mathcal{A}}^{(2)}}\rrbracket}$.

We next introduce the notion of second-order path term. For a sort $s\in S$, a term in $\mathrm{T}_{\Sigma^{\boldsymbol{\mathcal{A}}^{(2)}}}(X)_{s}$ will be called second-order path term if it is a term that can be interpreted as a  second-order path in $\mathrm{Pth}_{\boldsymbol{\mathcal{A}}^{(2)},s}$ by means of the $\mathrm{ip}^{(2,X)@}$ mapping. We will see that these terms play a crucial role in iterating the results of the first part.

\begin{restatable}{definition}{DDPT}
\label{DDPT} 
\index{path terms!second-order!$\mathrm{PT}_{\boldsymbol{\mathcal{A}}^{(2)}}$}
The $S$-sorted set of \emph{second-order path terms}, denoted by $\mathrm{PT}_{\boldsymbol{\mathcal{A}}^{(2)}}$, is the $S$-sorted subset of $\mathrm{T}_{\Sigma^{\boldsymbol{\mathcal{A}}^{(2)}}}(X)$ given by the $\Theta^{\llbracket 2 \rrbracket}$-saturation of $\mathrm{CH}^{(2)}[\mathrm{Pth}_{\boldsymbol{\mathcal{A}}^{(2)}}]$,
$$
\mathrm{PT}_{\boldsymbol{\mathcal{A}}^{(2)}}=\left[\mathrm{CH}^{(2)}[\mathrm{Pth}_{\boldsymbol{\mathcal{A}}^{(2)}}]\right]^{\Theta^{\llbracket 2 \rrbracket}}.
$$

Thus, for a sort $s\in S$, a term $P\in\mathrm{T}_{\Sigma^{\boldsymbol{\mathcal{A}}^{(2)}}}(X)_{s}$ belongs to  $\mathrm{PT}_{\boldsymbol{\mathcal{A}}^{(2)},s}$ if, and only if, there exists a second-order path $\mathfrak{P}^{(2)}\in\mathrm{Pth}_{\boldsymbol{\mathcal{A}}^{(2)},s}$ such that 
$$\left(
P,\mathrm{CH}^{(2)}_{s}\left(
\mathfrak{P}^{(2)}
\right)\right)
\in\Theta^{\llbracket 2 \rrbracket}_{s}.$$
\end{restatable}

In the following proposition we prove the most important property of second-order path terms. For every sort $s\in S$, a term $P$ in $\mathrm{T}_{\Sigma^{\boldsymbol{\mathcal{A}}^{(2)}}}(X)_{s}$ is a second-order path term if, and only if, $\mathrm{ip}^{(2,X)@}_{s}(P)$ is a second-order path.

\begin{restatable}{proposition}{PDPT}
\label{PDPT} Let $s$ be a sort in $S$ and $P$ a term in $\mathrm{T}_{\Sigma^{\boldsymbol{\mathcal{A}}^{(2)}}}(X)_{s}$, then the following conditions are equivalent.
\begin{enumerate}
\item[(i)] $P\in\mathrm{PT}_{\boldsymbol{\mathcal{A}}^{(2)},s}$;
\item[(ii)] $\mathrm{ip}^{(2,X)@}_{s}(P)$ is a second-order path in $\mathrm{Pth}_{\boldsymbol{\mathcal{A}}^{(2)},s}$.
\end{enumerate} 
\end{restatable}
\begin{proof}
Assume that $P$ is a term in $\mathrm{PT}_{\boldsymbol{\mathcal{A}}^{(2)},s}$. Then there exists a second-order path $\mathfrak{P}^{(2)}\in\mathrm{Pth}_{\boldsymbol{\mathcal{A}}^{(2)},s}$ for which $(P,\mathrm{CH}^{(2)}_{s}(\mathfrak{P}^{(2)}))$ is a pair in $\Theta^{\llbracket 2 \rrbracket}_{s}$. Then, by Corollary~\ref{CDVCong}, $\mathrm{ip}^{(2,X)@}_{s}(P)$ is a second-order path in $\mathrm{Pth}_{\boldsymbol{\mathcal{A}}^{(2)},s}$. 

Now, assume that $\mathrm{ip}^{(2,X)@}_{s}(P)$ is a second-order path in $\mathrm{Pth}_{\boldsymbol{\mathcal{A}}^{(2)},s}$. Then, by Lemma~\ref{LDVWCong}, $(P,\mathrm{CH}^{(2)}_{s}(\mathrm{ip}^{(2,X)@}_{s}(P)))$ is a pair in $\Theta^{\llbracket 2 \rrbracket}_{s}$, i.e., $P$ is a term in $\mathrm{PT}_{\boldsymbol{\mathcal{A}}^{(2)},s}$.
\end{proof}

In the following corollary we state that two second-order path terms in the same $\Theta^{\llbracket 2 \rrbracket}$-class, when interpreted as second-order paths by means of the  $\mathrm{ip}^{(2,X)@}$ mapping, have the same $([1],2)$-source and $([1],2)$-target and, consequently, the same $(0,2)$-source and $(0,2)$-target.

\begin{restatable}{corollary}{CDPTScTg}
\label{CDPTScTg}
Let $s$ be a sort in $S$ and $Q$, $P$ second-order path terms in $\mathrm{PT}_{\boldsymbol{\mathcal{A}}^{(2)},s}$ satisfying that
$(Q,P)\in\Theta^{\llbracket 2 \rrbracket}_{s}$. Then 
\begin{align*}
\left(
\mathrm{ip}^{(2,X)@}_{s}\left(
Q
\right),
\mathrm{ip}^{(2,X)@}_{s}\left(
P
\right)\right)
&\in\mathrm{Ker}\left(
\mathrm{sc}^{([1],2)}
\right)_{s}
\cap\mathrm{Ker}\left(
\mathrm{tg}^{([1],2)}
\right)_{s};
\\
\left(
\mathrm{ip}^{(2,X)@}_{s}\left(
Q\right),
\mathrm{ip}^{(2,X)@}_{s}\left(
P
\right)\right)
&\in\mathrm{Ker}\left(
\mathrm{sc}^{(0,2)}
\right)_{s}\cap\mathrm{Ker}\left(
\mathrm{tg}^{(0,2)}
\right)_{s}.
\end{align*}
\end{restatable}
\begin{proof}
By Proposition~\ref{PDPT} and Lemma~\ref{LDVCong}, it is the case that 
$$
\Bigl\llbracket
\mathrm{ip}^{(2,X)@}_{s}\left(
Q
\right)
\Bigr\rrbracket_{s}
=
\Bigl\llbracket
\mathrm{ip}^{(2,X)@}_{s}\left(
P
\right)
\Bigr\rrbracket_{s}.
$$

This implies, by Lemma~\ref{LDV}, that 
$$
\left(
\mathrm{ip}^{(2,X)@}_{s}\left(
Q\right),
\mathrm{ip}^{(2,X)@}_{s}\left(
P
\right)\right)
\in\mathrm{Ker}\left(
\mathrm{sc}^{([1],2)}
\right)_{s}\cap\mathrm{Ker}\left(
\mathrm{tg}^{([1],2)}
\right)_{s}.
$$

From the last item and taking into account the definition of the $(0,2)$-source and $(0,2)$-target introduced in Definition~\ref{DDScTgZ} we conclude that 
$$
\left(
\mathrm{ip}^{(2,X)@}_{s}\left(
Q
\right),
\mathrm{ip}^{(2,X)@}_{s}\left(
P
\right)\right)
\in\mathrm{Ker}\left(
\mathrm{sc}^{(0,2)}
\right)_{s}
\cap\mathrm{Ker}\left(
\mathrm{tg}^{(0,2)}
\right)_{s}.
$$

This completes the proof.
\end{proof}

In the following corollaries we state that some already known $S$-sorted mappings have nice restrictions or correstrictions to the many-sorted set of path terms introduced in this chapter. We will start with the embeddings introduced in Definition~\ref{DDEta}.

\begin{restatable}{corollary}{CDPTEta}
\label{CDPTEta} The $S$-sorted embeddings
$\eta^{(2,X)}$, $\eta^{(2,\mathcal{A})}$ and $\eta^{(2,\mathcal{A}^{(2)})}$  introduced in Definition~\ref{DDEta} from, respectively, $X$, $\mathcal{A}$ and $\mathcal{A}^{(2)}$ to $\mathrm{T}_{\Sigma^{\boldsymbol{\mathcal{A}}^{(2)}}}(X)$ correstrict to $\mathrm{PT}_{\boldsymbol{\mathcal{A}}^{(2)}}$.
\end{restatable}
\begin{proof}
The following equations hold by Propositions~\ref{PDCHDUId},~\ref{PDCHA} and~\ref{PDCHDA}, respectively
\begin{enumerate}
\item $\eta^{(2,X)}=
\mathrm{CH}^{(2)}\circ\mathrm{ip}^{(2,X)};$
\item $
\eta^{(2,\mathcal{A})}=
\mathrm{CH}^{(2)}\circ\mathrm{ech}^{(2,\mathcal{A})};$
\item $\eta^{(2,\mathcal{A}^{(2)})}=
\mathrm{CH}^{(2)}\circ\mathrm{ech}^{(2,\mathcal{A}^{(2)})}.$
\end{enumerate}

Therefore, $\eta^{(2,X)}$, $\eta^{(2,\mathcal{A})}$ and $\eta^{(2,\mathcal{A}^{(2)})}$ always retrieve terms in $\mathrm{CH}^{(2)}[\mathrm{Pth}_{\boldsymbol{\mathcal{A}}^{(2)}}]$.
\end{proof}

According to the above corollary, we can consider the correstriction of the embeddings introduced in Definition~\ref{DDEta} to second-order path terms. To avoid further notation, we resignify the already existing many-sorted mappings.

\begin{restatable}{definition}{DDPTEta}
\label{DDPTEta}  Let $X$ be an $S$-sorted set and let $\mathrm{PT}_{\boldsymbol{\mathcal{A}}^{(2)}}$ be the many-sorted set of second-order path terms.  From now on, we will denote by

\begin{enumerate}
\item $\eta^{(2,X)}$ the $S$-sorted mapping from $X$ to 
$\mathrm{PT}_{\boldsymbol{\mathcal{A}}^{(2)}}$ such that, for every sort $s\in S$, sends an element $x\in X_{s}$ to the variable  $\eta^{(2,X)}_{s}(x)$ in $\mathrm{PT}_{\boldsymbol{\mathcal{A}}^{(2)},s}$.
\index{inclusion!second-order!$\eta^{(2,X)}$}
\item $\eta^{(2,\mathcal{A})}$ the $S$-sorted mapping from 
$\mathcal{A}$ to 
$\mathrm{PT}_{\boldsymbol{\mathcal{A}}^{(2)}}$ 
such that, for every sort $s\in S$, sends a  rewrite rule $\mathfrak{p}\in \mathcal{A}_{s}$ to the constant $\mathfrak{p}^{\mathrm{T}_{\Sigma^{\boldsymbol{\mathcal{A}}^{(2)}}}(X)}$.
\index{inclusion!second-order!$\eta^{(2,\mathcal{A})}$}
\item $\eta^{(2,\mathcal{A}^{(2)})}$ the $S$-sorted mapping from 
$\mathcal{A}^{(2)}$ to 
$\mathrm{PT}_{\boldsymbol{\mathcal{A}}^{(2)}}$ 
such that, for every sort $s\in S$, sends a second-order rewrite rule $\mathfrak{p}^{(2)}\in \mathcal{A}^{(2)}_{s}$ to the constant $\mathfrak{p}^{(2)\mathrm{T}_{\Sigma^{\boldsymbol{\mathcal{A}}^{(2)}}}(X)}$.
\index{inclusion!second-order!$\eta^{(2,\mathcal{A}^{(2)})}$}
\end{enumerate} 

The above $S$-sorted mappings are depicted in the diagram of Figure~\ref{FDPTEta}.
\end{restatable}

\begin{figure}
\begin{tikzpicture}
[ACliment/.style={-{To [angle'=45, length=5.75pt, width=4pt, round]}},scale=.8]
\node[] (x) at (0,0) [] {$X$};
\node[] (a) at (0,-1.5) [] {$\mathcal{A}$};
\node[] (a2) at (0,-3) [] {$\mathcal{A}^{(2)}$};
\node[] (T) at (6,-3) [] {$\mathrm{PT}_{\boldsymbol{\mathcal{A}}^{(2)}}$};
\draw[ACliment, bend left=20]  (x) to node [above right] {$\eta^{(2,X)}$} (T);
\draw[ACliment, bend left=10]  (a) to node [midway, fill=white] {$\eta^{(2,\mathcal{A})}$} (T);
\draw[ACliment]  (a2) to node [below] {$\eta^{(2,\mathcal{A}^{(2)})}$} (T);
\end{tikzpicture}
\caption{Refined embeddings relative to $X$, $\mathcal{A}$ and $\mathcal{A}^{(2)}$ at layer 2.}\label{FDPTEta}
\end{figure}

We next consider  the embeddings introduced in Propositions~\ref{PDUEmb} and~\ref{PDZEmb}.

\begin{restatable}{corollary}{CDPTEtas}
\label{CDPTEtas} The $S$-sorted embeddings $\eta^{(2,1)\sharp}$ and $\eta^{(2,0)\sharp}$ introduced in Propositions~\ref{PDUEmb} and~\ref{PDZEmb} from, respectively, $\mathrm{PT}_{\boldsymbol{\mathcal{A}}}$ and $\mathrm{T}_{\Sigma}(X)$ to $\mathrm{T}_{\Sigma^{\boldsymbol{\mathcal{A}}^{(2)}}}(X)$ correstrict to $\mathrm{PT}_{\boldsymbol{\mathcal{A}}^{(2)}}$.
\end{restatable}
\begin{proof}
The following equalities hold by, respectively, Propositions~\ref{PDIpDU} and~\ref{PDIpDZ}
\begin{enumerate}
\item $\mathrm{ip}^{(2,X)@}\circ\eta^{(2,1)\sharp}=
\eta^{(2,\mathrm{Pth}_{\boldsymbol{\mathcal{A}}^{(2)}})}\circ\mathrm{ip}^{(2,[1])\sharp}\circ\mathrm{pr}^{\Theta^{[1]}}
$;
\item $\mathrm{ip}^{(2,X)@}\circ\eta^{(2,0)\sharp}=
\eta^{(2,\mathrm{Pth}_{\boldsymbol{\mathcal{A}}^{(2)}})}\circ\mathrm{ip}^{(2,0)\sharp}$.
\end{enumerate}

According to (1), for every sort $s\in S$, and every path term $P\in\mathrm{PT}_{\boldsymbol{\mathcal{A}},s}$, we have that $\mathrm{ip}^{(2,X)@}_{s}(\eta^{(2,1)\sharp}_{s}(P))$ is a second-order path in $\mathrm{Pth}_{\boldsymbol{\mathcal{A}}^{(2)},s}$. We conclude following Proposition~\ref{PDPT} that $\eta^{(2,1)\sharp}_{s}(P)$ is a second-order path term in $\mathrm{PT}_{\boldsymbol{\mathcal{A}}^{(2)},s}$.

According to (2), for every sort $s\in S$, and every term $P\in\mathrm{T}_{\Sigma}(X)$, we have that $\mathrm{ip}^{(2,X)@}_{s}(\eta^{(2,0)\sharp}_{s}(P))$ is a second-order path in $\mathrm{Pth}_{\boldsymbol{\mathcal{A}}^{(2)},s}$. We conclude following Proposition~\ref{PDPT} that $\eta^{(2,0)\sharp}_{s}(P)$ is a second-order path term in $\mathrm{PT}_{\boldsymbol{\mathcal{A}}^{(2)},s}$.

This completes the proof.
\end{proof}

According to the above corollary, we can consider the correstriction of the embeddings introduced in Propositions~\ref{PDUEmb} and~\ref{PDZEmb} to second-order path terms. To avoid further notation, we resignify the already existing many-sorted mappings.

\begin{restatable}{definition}{DDPTEtas}
\label{DDPTEtas}  Let $X$ be an $S$-sorted set and let $\mathrm{PT}_{\boldsymbol{\mathcal{A}}^{(2)}}$ be the many-sorted set of second-order path terms.  From now on, we will denote by

\begin{enumerate}
\item  $\eta^{(2,1)\sharp}$ the $S$-sorted mapping from $\mathrm{PT}_{\boldsymbol{\mathcal{A}}}$ to 
$\mathrm{PT}_{\boldsymbol{\mathcal{A}}^{(2)}}$ such that, for every sort $s\in S$, sends a path term $P\in \mathrm{PT}_{\boldsymbol{\mathcal{A}},s}$ to the second-order path term  $
\eta^{(2,1)\sharp}_{s}(P)$ in $\mathrm{PT}_{\boldsymbol{\mathcal{A}}^{(2)},s}$.
\index{inclusion!second-order!$\eta^{(2,1)\sharp}$}
\item $\eta^{(2,0)\sharp}$ the $S$-sorted mapping from $\mathrm{T}_{\Sigma}(X)$ to 
$\mathrm{PT}_{\boldsymbol{\mathcal{A}}^{(2)}}$ such that, for every sort $s\in S$, sends a term $P\in \mathrm{T}_{\Sigma}(X)_{s}$ to the second-order path term  $\eta^{(2,0)\sharp}_{s}(P)$ in $\mathrm{PT}_{\boldsymbol{\mathcal{A}}^{(2)},s}$.
\index{inclusion!second-order!$\eta^{(2,0)\sharp}$}
\end{enumerate}
\end{restatable}

We next consider the case of the $\mathrm{ip}^{(2,X)@}$ mapping.

\begin{restatable}{corollary}{CDPTIp}
\label{CDPTIp} 
The restriction of the  $\mathrm{ip}^{(2,X)@}$ mapping to the set of second-order path terms, that is,
$$
\mathrm{ip}^{(2,X)@}{\,\!\upharpoonright}_{\mathrm{PT}_{\boldsymbol{\mathcal{A}}^{(2)}}}
\colon
\mathrm{PT}_{\boldsymbol{\mathcal{A}}^{(2)}}
\mor
\mathrm{F}_{\Sigma^{\boldsymbol{\mathcal{A}}^{(2)}}}(\mathbf{Pth}_{\boldsymbol{\mathcal{A}}^{(2)}}),
$$
correstricts to $\mathrm{Pth}_{\boldsymbol{\mathcal{A}}^{(2)}}$.
\end{restatable}
\begin{proof}
It follows from Proposition~\ref{PDPT}.
\end{proof}

\begin{definition}
\index{identity!second-order!$\mathrm{ip}^{(2,X)a}$}
  From now on, we resignify  $\mathrm{ip}^{(2,X)@}$ to match the birrestriction introduced in Corollary~\ref{CDPTIp}. Thus,
$$
\mathrm{ip}^{(2,X)@}
\colon
\mathrm{PT}_{\boldsymbol{\mathcal{A}}^{(2)}}
\mor
\mathrm{Pth}_{\boldsymbol{\mathcal{A}}^{(2)}}.
$$
\end{definition}

Now it is the turn of the second-order Curry-Howard mapping.

\begin{restatable}{corollary}{CDPTCH}
\label{CDPTCH}
The second-order Curry-Howard mapping, that is,
$$
\mathrm{CH}^{(2)}\colon\mathrm{Pth}_{\boldsymbol{\mathcal{A}}^{(2)}}
\mor
\mathrm{T}_{\Sigma^{\boldsymbol{\mathcal{A}}^{(2)}}}(X),
$$
correstricts to $\mathrm{PT}_{\boldsymbol{\mathcal{A}}^{(2)}}$.
\end{restatable}
\begin{proof}
Note that $\mathrm{CH}^{(2)}[\mathrm{Pth}_{\boldsymbol{\mathcal{A}}^{(2)}}]$ is a subset of its own saturation with respect to the $\Sigma^{\boldsymbol{\mathcal{A}}^{(2)}}$-congruence $\Theta^{\llbracket 2 \rrbracket}$.
\end{proof}

\begin{definition}
\index{Curry-Howard!second-order!$\mathrm{CH}^{(2)}$}
 From now on, we resignify  $\mathrm{CH}^{(2)}$ to match the  correstriction presented in Corollary~\ref{CDPTCH}, that is,
$$
\mathrm{CH}^{(2)}\colon\mathrm{Pth}_{\boldsymbol{\mathcal{A}}^{(2)}}
\mor
\mathrm{PT}_{\boldsymbol{\mathcal{A}}^{(2)}}.
$$

\end{definition}

\begin{figure}
\begin{tikzpicture}
[ACliment/.style={-{To [angle'=45, length=5.75pt, width=4pt, round]}},scale=0.8]
\node[] (x) at (0,0) [] {$X$};
\node[] (t) at (6,0) [] {$\mathrm{T}_{\Sigma}(X)$};
\node[] (pth1) at (6,-2) [] {$\mathrm{Pth}_{\boldsymbol{\mathcal{A}}}$};
\node[] (pt1) at (6,-4) [] {$\mathrm{PT}_{\boldsymbol{\mathcal{A}}}$};

\node[] (pth2) at (6,-6) [] {$\mathrm{Pth}_{\boldsymbol{\mathcal{A}}^{(2)}}$};
\node[] (pt2) at (6,-8) [] {$\mathrm{PT}_{\boldsymbol{\mathcal{A}}^{(2)}}$};

\draw[ACliment]  (x) to node [above right]
{$\textstyle \eta^{(0,X)}$} (t);
\draw[ACliment, bend right=10]  (x) to node [pos=.7, fill=white]
{$\textstyle \mathrm{ip}^{(1,X)}$} (pth1.west);
\draw[ACliment, bend right=20]  (x) to node  [pos=.6, fill=white]
{$\textstyle\eta^{(1,X)}$} (pt1.west);

\draw[ACliment, bend right=30]  (x) to node [pos=.5, fill=white]
{$\textstyle \mathrm{ip}^{(2,X)}$} (pth2.west);
\draw[ACliment, bend right=40]  (x) to node  [below left]
{$\textstyle \eta^{(2,X)}$} (pt2.west);

\draw[ACliment] 
($(t)+(0,-.35)$) to node [right] {$\textstyle \mathrm{ip}^{(1,0)\sharp}$}  ($(pth1)+(0,.35)$);

\draw[ACliment] 
($(pth1)+(-.30,-.35)$) to node [left] {$\textstyle \mathrm{CH}^{(1)}$}  ($(pt1)+(-.30,.35)$);
\draw[ACliment] 
($(pt1)+(.30,+.35)$) to node [right] {$\textstyle \mathrm{ip}^{(1,X)@}$}  ($(pth1)+(.30,-.35)$);


\draw[ACliment] 
($(pt1)+(0,-.35)$) to node [right] {$\textstyle \mathrm{ip}^{(2,1)\sharp}$}  ($(pth2)+(0,.35)$);

\draw[ACliment] 
($(pth2)+(-.30,-.35)$) to node [left] {$\textstyle \mathrm{CH}^{(2)}$}  ($(pt2)+(-.30,.35)$);
\draw[ACliment] 
($(pt2)+(.30,+.35)$) to node [right] {$\textstyle \mathrm{ip}^{(2,X)@}$}  ($(pth2)+(.30,-.35)$);

\draw[ACliment, rounded corners]
(pt1.east) -- 
($(pt1)+(2,0)$) -- node [right] {$\textstyle \eta^{(2,1)\sharp}$} 
($(pt2)+(2,+.1)$) -- ($(pt2.east)+(0,.1)$)
;

\draw[ACliment, rounded corners]
(t.east) -- 
($(t)+(3.5,0)$) -- node [right] {$\textstyle \eta^{(2,0)\sharp}$} 
($(pt2)+(3.5,-.1)$) -- 
($(pt2.east)-(0,.1)$)
;

\end{tikzpicture}
\caption{Refined embeddings relative to $X$ at layers $0$, $1$ \& $2$.}
\label{FDPTX}
\end{figure}

\begin{figure}
\begin{tikzpicture}
[ACliment/.style={-{To [angle'=45, length=5.75pt, width=4pt, round]}},scale=0.8]
\node[] (a) at (0,-2) [] {$\mathcal{A}$};
\node[] (pth1) at (6,-2) [] {$\mathrm{Pth}_{\boldsymbol{\mathcal{A}}}$};
\node[] (pt1) at (6,-4) [] {$\mathrm{PT}_{\boldsymbol{\mathcal{A}}}$};

\node[] (pth2) at (6,-6) [] {$\mathrm{Pth}_{\boldsymbol{\mathcal{A}}^{(2)}}$};
\node[] (pt2) at (6,-8) [] {$\mathrm{PT}_{\boldsymbol{\mathcal{A}}^{(2)}}$};

\draw[ACliment]  (a) to node [above right]
{$\textstyle \mathrm{ech}^{(1,\mathcal{A})}$} (pth1);
\draw[ACliment, bend right=10]  (a) to node [pos=.6, fill=white]
{$\textstyle \eta^{(1,\mathcal{A})}$} (pt1.west);
\draw[ACliment, bend right=20]  (a) to node  [midway, fill=white]
{$\textstyle \mathrm{ech}^{(2,\mathcal{A})}$} (pth2.west);

\draw[ACliment, bend right=30]  (a) to node [below left]
{$\textstyle \eta^{(2,\mathcal{A})}$} (pt2.west);

\draw[ACliment] 
($(pth1)+(-.30,-.35)$) to node [left] {$\textstyle \mathrm{CH}^{(1)}$}  ($(pt1)+(-.30,.35)$);
\draw[ACliment] 
($(pt1)+(.30,+.35)$) to node [right] {$\textstyle \mathrm{ip}^{(1,X)@}$}  ($(pth1)+(.30,-.35)$);


\draw[ACliment] 
($(pt1)+(0,-.35)$) to node [right] {$\textstyle \mathrm{ip}^{(2,1)\sharp}$}  ($(pth2)+(0,.35)$);

\draw[ACliment] 
($(pth2)+(-.30,-.35)$) to node [left] {$\textstyle \mathrm{CH}^{(2)}$}  ($(pt2)+(-.30,.35)$);
\draw[ACliment] 
($(pt2)+(.30,+.35)$) to node [right] {$\textstyle \mathrm{ip}^{(2,X)@}$}  ($(pth2)+(.30,-.35)$);

\draw[ACliment, rounded corners]
(pt1.east) -- 
($(pt1)+(2,0)$) -- node [right] {$\textstyle \eta^{(2,1)\sharp}$} 
($(pt2)+(2,0)$) -- ($(pt2.east)+(0,0)$)
;

\end{tikzpicture}
\caption{Refined embeddings relative to $\mathcal{A}$ at layers $1$ \& $2$.}
\label{FDPTA}
\end{figure}

\begin{figure}
\begin{tikzpicture}
[ACliment/.style={-{To [angle'=45, length=5.75pt, width=4pt, round]}},scale=0.8]
\node[] (a) at (0,-2) [] {$\mathcal{A}^{(2)}$};
\node[] (pth1) at (6,-2) [] {$\mathrm{Pth}_{\boldsymbol{\mathcal{A}}^{(2)}}$};
\node[] (pt1) at (6,-4) [] {$\mathrm{PT}_{\boldsymbol{\mathcal{A}}^{(2)}}$};

\draw[ACliment]  (a) to node [above]
{$\textstyle \mathrm{ech}^{(2,\mathcal{A}^{(2)})}$} (pth1);
\draw[ACliment, bend right=10]  (a) to node  [below left]
{$\textstyle \eta^{(2,\mathcal{A}^{(2)})}$} (pt1.west);

\draw[ACliment] 
($(pth1)+(-.3,-.35)$) to node [left] {$\textstyle \mathrm{CH}^{(2)}$}  ($(pt1)+(-.3,.35)$);
\draw[ACliment] 
($(pt1)+(.3,+.35)$) to node [right] {$\textstyle \mathrm{ip}^{(2,X)@}$}  ($(pth1)+(.3,-.35)$);
\end{tikzpicture}
\caption{Refined embeddings relative to $\mathcal{A}^{(2)}$ at layer $2$.}
\label{FDPTDA}
\end{figure}

\begin{remark} The above refinements are considered in Figures~\ref{FDPTX},~\ref{FDPTA} and~\ref{FDPTDA}, relatively to $X$, $\mathcal{A}$ and $\mathcal{A}^{(2)}$, respectively. Every possible triangle considered in these figures commute.
\end{remark}

\section{
\texorpdfstring
{A structure of partial $\Sigma^{\boldsymbol{\mathcal{A}}^{(2)}}$-algebra on $\mathrm{PT}_{\boldsymbol{\mathcal{A}}^{(2)}}$}
{A partial algebra on second-order path terms}
}

We next show that the many-sorted set of second-order path terms has a natural structure of partial $\Sigma^{\boldsymbol{\mathcal{A}}^{(2)}}$-algebra, which is a subalgebra of $\mathbf{T}_{\Sigma^{\boldsymbol{\mathcal{A}}^{(2)}}}(X)$. 

\begin{restatable}{proposition}{PDPTDCatAlg}
\label{PDPTDCatAlg} 
\index{path terms!second-order!$\mathbf{PT}_{\boldsymbol{\mathcal{A}}^{(2)}}$}
The $S$-sorted set $\mathrm{PT}_{\boldsymbol{\mathcal{A}}^{(2)}}$ is equipped, in a natural way, with a structure of many-sorted partial $\Sigma^{\boldsymbol{\mathcal{A}}^{(2)}}$-algebra, which is a $\Sigma^{\boldsymbol{\mathcal{A}}^{(2)}}$-subalgebra of $\mathbf{T}_{\Sigma^{\boldsymbol{\mathcal{A}}^{(2)}}}(X)$.
\end{restatable}
\begin{proof}
Let us denote by $\mathbf{PT}_{\boldsymbol{\mathcal{A}}^{(2)}}$ the partial $\Sigma^{\boldsymbol{\mathcal{A}}^{(2)}}$-algebra defined as follows:

\textsf{(i)} 
The underlying $S$-sorted set of $\mathbf{PT}_{\boldsymbol{\mathcal{A}}^{(2)}}$ is $\mathrm{PT}_{\boldsymbol{\mathcal{A}}^{(2)}}$.

\textsf{(ii)} For every $(\mathbf{s},s)\in S^{\star}\times S$ and every operation symbol $\sigma\in\Sigma_{\mathbf{s},s}$, the operation $\sigma^{\mathbf{PT}_{\boldsymbol{\mathcal{A}}^{(2)}}}$ is equal to $\sigma^{\mathbf{T}_{\Sigma^{\boldsymbol{\mathcal{A}}^{(2)}}}(X)}$, i.e., to the interpretation of $\sigma$ in $\mathbf{T}_{\Sigma^{\boldsymbol{\mathcal{A}}^{(2)}}}(X)$. 

Let us note that, if $(P_{j})_{j\in\bb{\mathbf{s}}}$ is a family of second-order path terms in $\mathrm{PT}_{\boldsymbol{\mathcal{A}}^{(2)},\mathbf{s}}$, then the following chain of equalities holds
\allowdisplaybreaks
\begin{align*}
\mathrm{ip}^{(2,X)@}_{s}\left(
\sigma^{\mathbf{T}_{\Sigma^{\boldsymbol{\mathcal{A}}^{(2)}}}(X)}
\left(\left(
P_{j}
\right)_{j\in\bb{\mathbf{s}}}
\right)\right)
&=
\sigma^{\mathbf{F}_{\Sigma^{\boldsymbol{\mathcal{A}}^{(2)}}}\left(
\mathbf{Pth}_{\boldsymbol{\mathcal{A}}^{(2)}}
\right)}
\left(\left(
\mathrm{ip}^{(2,X)@}_{s_{j}}\left(
P_{j}
\right)\right)_{j\in\bb{\mathbf{s}}}
\right)
\tag{1}
\\&=
\sigma^{\mathbf{Pth}_{\boldsymbol{\mathcal{A}}^{(2)}}}
\left(\left(
\mathrm{ip}^{(2,X)@}_{s_{j}}\left(
P_{j}
\right)\right)_{j\in\bb{\mathbf{s}}}\right).
\tag{2}
\end{align*}

The first equality follows from the fact that, by Definition~\ref{DDIp}, $\mathrm{ip}^{(2,X)@}$ is a $\Sigma^{\boldsymbol{\mathcal{A}}^{(2)}}$-homomorphism; finally, the second equality follows from the fact that, for every $j\in\bb{\mathbf{s}}$, the term $P_{j}$ is a second-order path term in $\mathrm{PT}_{\boldsymbol{\mathcal{A}}^{(2)},s_{j}}$. Hence, in virtue of Proposition~\ref{PDPT}, for every $j\in\bb{\mathbf{s}}$, $\mathrm{ip}^{(2,X)@}_{s_{j}}(P_{j})$ is a second-order path in $\mathrm{Pth}_{\boldsymbol{\mathcal{A}}^{(2)},s_{j}}$. Thus, the interpretation of the operation symbol $\sigma$ in $\mathbf{F}_{\Sigma^{\boldsymbol{\mathcal{A}}^{(2)}}}(\mathbf{Pth}_{\boldsymbol{\mathcal{A}}^{(2)}})$ is given by the interpretation of $\sigma$ in $\mathbf{Pth}_{\boldsymbol{\mathcal{A}}^{(2)}}$.

Hence, since $\mathrm{ip}^{(2,X)@}_{s}(\sigma^{\mathbf{T}_{\Sigma^{\boldsymbol{\mathcal{A}}^{(2)}}}(X)}((P_{j})_{j\in\bb{\mathbf{s}}}))$ is a second-order path in $\mathrm{Pth}_{\boldsymbol{\mathcal{A}}^{(2)},s}$, we conclude, following Proposition~\ref{PDPT}, that $\sigma^{\mathbf{T}_{\Sigma^{\boldsymbol{\mathcal{A}}^{(2)}}}(X)}((P_{j})_{j\in\bb{\mathbf{s}}})$ is a second-order path term in $\mathrm{PT}_{\boldsymbol{\mathcal{A}}^{(2)},s}$. That is, the operation $\sigma^{\mathbf{PT}_{\boldsymbol{\mathcal{A}}^{(2)}}}$ is well-defined.

\textsf{(iii)} For every $s\in S$ and every rewrite rule $\mathfrak{p}\in \mathcal{A}_{s}$, the constant $\mathfrak{p}^{\mathbf{PT}_{\boldsymbol{\mathcal{A}}^{(2)}}}$ is equal to $\mathfrak{p}^{\mathbf{T}_{\Sigma^{\boldsymbol{\mathcal{A}}^{(2)}}}(X)}$, i.e., to the interpretation of $\mathfrak{p}$ in $\mathbf{T}_{\Sigma^{\boldsymbol{\mathcal{A}}^{(2)}}}(X)$. 

Let us note that the following chain of equalities holds
\allowdisplaybreaks
\begin{align*}
\mathrm{ip}^{(2,X)@}_{s}\left(
\mathfrak{p}^{\mathbf{T}_{\Sigma^{\boldsymbol{\mathcal{A}}^{(2)}}}(X)}
\right)
&=
\mathfrak{p}^{\mathbf{F}_{\Sigma^{\boldsymbol{\mathcal{A}}^{(2)}}}\left(
\mathbf{Pth}_{\boldsymbol{\mathcal{A}}^{(2)}}\right)}
\tag{1}
\\&=
\mathfrak{p}^{\mathbf{Pth}_{\boldsymbol{\mathcal{A}}^{(2)}}}.
\tag{2}
\end{align*}

The first equality follows from the fact that, by Definition~\ref{DDIp}, $\mathrm{ip}^{(2,X)@}$ is a $\Sigma^{\boldsymbol{\mathcal{A}}^{(2)}}$-homomorphism; finally, the second equality follows from the fact that the interpretation of the constant operation symbol $\mathfrak{p}$ in $\mathbf{F}_{\Sigma^{\boldsymbol{\mathcal{A}}^{(2)}}}(\mathbf{Pth}_{\boldsymbol{\mathcal{A}}^{(2)}})$ is given by the interpretation of $\mathfrak{p}$ in $\mathbf{Pth}_{\boldsymbol{\mathcal{A}}^{(2)}}$.

Thus, since $\mathrm{ip}^{(2,X)@}_{s}(\mathfrak{p}^{\mathbf{T}_{\Sigma^{\boldsymbol{\mathcal{A}}^{(2)}}}(X)})$ is a second-order path in $\mathrm{Pth}_{\boldsymbol{\mathcal{A}}^{(2)},s}$, we conclude, by Proposition~\ref{PDPT}, that $\mathfrak{p}^{\mathbf{T}_{\Sigma^{\boldsymbol{\mathcal{A}}^{(2)}}}(X)}$ is a second-order path term in $\mathrm{PT}_{\boldsymbol{\mathcal{A}}^{(2)},s}$. That is, the constant $\mathfrak{p}^{\mathbf{T}_{\Sigma^{\boldsymbol{\mathcal{A}}^{(2)}}}(X)}$ is well-defined.

\textsf{(iv)} For every $s\in S$, the operation $\mathrm{sc}_{s}^{0\mathbf{PT}_{\boldsymbol{\mathcal{A}}^{(2)}}}$ is equal to $\mathrm{sc}_{s}^{0\mathbf{T}_{\Sigma^{\boldsymbol{\mathcal{A}}^{(2)}}}(X)}$, i.e., to the interpretation of $\mathrm{sc}_{s}^{0}$ in $\mathbf{T}_{\Sigma^{\boldsymbol{\mathcal{A}}^{(2)}}}(X)$. 

Let us note that, if $P$ is a second-order path term in $\mathrm{PT}_{\boldsymbol{\mathcal{A}}^{(2)},s}$, then the following chain of equalities holds
\allowdisplaybreaks
\begin{align*}
\mathrm{ip}^{(2,X)@}_{s}\left(
\mathrm{sc}_{s}^{0\mathbf{T}_{\Sigma^{\boldsymbol{\mathcal{A}}^{(2)}}}(X)}\left(
P
\right)\right)
&=
\mathrm{sc}_{s}^{0\mathbf{F}_{\Sigma^{\boldsymbol{\mathcal{A}}^{(2)}}}\left(
\mathbf{Pth}_{\boldsymbol{\mathcal{A}}^{(2)}}\right)}
\left(
\mathrm{ip}^{(2,X)@}_{s}\left(
P
\right)\right)
\tag{1}
\\&=
\mathrm{sc}_{s}^{0\mathbf{Pth}_{\boldsymbol{\mathcal{A}}^{(2)}}}\left(
\mathrm{ip}^{(2,X)@}_{s}\left(
P
\right)\right).
\tag{2}
\end{align*}

The first equality follows from the fact that, by Definition~\ref{DDIp}, $\mathrm{ip}^{(2,X)@}$ is a $\Sigma^{\boldsymbol{\mathcal{A}}^{(2)}}$-homomorphism;  finally, the second equality follows from the fact that $P$ is a second-order path term in $\mathrm{PT}_{\boldsymbol{\mathcal{A}}^{(2)},s}$. Hence, in virtue of Proposition~\ref{PDPT} $\mathrm{ip}^{(2,X)@}_{s}(P)$ is a second-order path in $\mathrm{Pth}_{\boldsymbol{\mathcal{A}}^{(2)},s}$. Thus, the interpretation of the operation symbol $\mathrm{sc}^{0}_{s}$ in $\mathbf{F}_{\Sigma^{\boldsymbol{\mathcal{A}}^{(2)}}}(\mathbf{Pth}_{\boldsymbol{\mathcal{A}}^{(2)}})$ is given by the interpretation of $\mathrm{sc}^{0}_{s}$ in $\mathbf{Pth}_{\boldsymbol{\mathcal{A}}^{(2)}}$.

Therefore, since $\mathrm{ip}^{(2,X)@}_{s}(\mathrm{sc}_{s}^{0\mathbf{T}_{\Sigma^{\boldsymbol{\mathcal{A}}^{(2)}}}(X)}(P))$ is a second-order path in $\mathrm{Pth}_{\boldsymbol{\mathcal{A}}^{(2)},s}$, we have that, by Proposition~\ref{PDPT}, $\mathrm{sc}_{s}^{0\mathbf{T}_{\Sigma^{\boldsymbol{\mathcal{A}}^{(2)}}}(X)}(P)$ is a second-order path term in $\mathrm{PT}_{\boldsymbol{\mathcal{A}}^{(2)},s}$. That is, the operation $\mathrm{sc}_{s}^{0\mathbf{PT}_{\boldsymbol{\mathcal{A}}^{(2)}}}$ is well-defined.

\textsf{(v)} For every $s\in S$, the operation $\mathrm{tg}_{s}^{0\mathbf{PT}_{\boldsymbol{\mathcal{A}}^{(2)}}}$ is equal to $\mathrm{tg}_{s}^{0\mathbf{T}_{\Sigma^{\boldsymbol{\mathcal{A}}^{(2)}}}(X)}$, i.e., to the interpretation of $\mathrm{tg}_{s}^{0}$ in $\mathbf{T}_{\Sigma^{\boldsymbol{\mathcal{A}}^{(2)}}}(X)$.  The operation $\mathrm{tg}_{s}^{0\mathbf{PT}_{\boldsymbol{\mathcal{A}}^{(2)}}}$ is well-defined by a similar argument to that of $\mathrm{sc}_{s}^{0\mathbf{PT}_{\boldsymbol{\mathcal{A}}^{(2)}}}$.

\textsf{(vi)} For every $s\in S$, the partial binary operation $\circ_{s}^{0\mathbf{PT}_{\boldsymbol{\mathcal{A}}^{(2)}}}$ is defined on second-order path terms $Q,P\in\mathrm{PT}_{\boldsymbol{\mathcal{A}}^{(2)},s}$ if, and only if, the following equation holds
$$
\mathrm{sc}_{s}^{(0,2)}\left(
\mathrm{ip}^{(2,X)@}_{s}\left(
Q
\right)\right)=
\mathrm{tg}_{s}^{(0,2)}\left(
\mathrm{ip}^{(2,X)@}_{s}\left(
P
\right)\right).
$$
For this case, $\circ_{s}^{0\mathbf{PT}_{\boldsymbol{\mathcal{A}}^{(2)}}}$ is equal to $\circ_{s}^{\mathbf{T}_{\Sigma^{\boldsymbol{\mathcal{A}}^{(2)}}}(X)}$, i.e., to the interpretation of $\circ^{0}_{s}$ in $\mathbf{T}_{\Sigma^{\boldsymbol{\mathcal{A}}^{(2)}}}(X)$.

Let us note that, if $Q$ and $P$ are second-order path terms in $\mathrm{PT}_{\boldsymbol{\mathcal{A}}^{(2)},s}$ with 
$$\mathrm{sc}_{s}^{(0,2)}\left(
\mathrm{ip}^{(2,X)@}_{s}\left(
Q
\right)\right)=
\mathrm{tg}_{s}^{(0,2)}\left(
\mathrm{ip}^{(2,X)@}_{s}\left(
P
\right)\right),$$
then the following chain of equalities holds
\begin{flushleft}
$\mathrm{ip}^{(2,X)@}_{s}\left(
Q\circ_{s}^{0\mathbf{T}_{\Sigma^{\boldsymbol{\mathcal{A}}^{(2)}}}(X)} P
\right)$
\allowdisplaybreaks
\begin{align*}
\qquad
&=
\left(
\mathrm{ip}^{(2,X)@}_{s}\left(
Q
\right)\right)
\circ_{s}^{0\mathbf{F}_{\Sigma^{\boldsymbol{\mathcal{A}}^{(2)}}}\left(
\mathbf{Pth}_{\boldsymbol{\mathcal{A}}^{(2)}}\right)}
\left(\mathrm{ip}^{(2,X)@}_{s}\left(
P
\right)\right)
\tag{1}
\\&=
\left(\mathrm{ip}^{(2,X)@}_{s}\left(Q
\right)\right)
\circ_{s}^{0\mathbf{Pth}_{\boldsymbol{\mathcal{A}}^{(2)}}}
\left(\mathrm{ip}^{(2,X)@}_{s}\left(
P
\right)\right).
\tag{2}
\end{align*}
\end{flushleft}

The first equality follows from the fact that, by Definition~\ref{DDIp}, $\mathrm{ip}^{(2,X)@}$ is a $\Sigma^{\boldsymbol{\mathcal{A}}^{(2)}}$-homomorphism;  finally, the second equality follows from the fact that $Q$ and $P$ are second-order path terms in $\mathrm{PT}_{\boldsymbol{\mathcal{A}}^{(2)},s}$. Hence, by Proposition~\ref{PDPT}, $\mathrm{ip}^{(2,X)@}_{s}(Q)$ and $\mathrm{ip}^{(2,X)@}_{s}(P)$ are second-order paths in $\mathrm{Pth}_{\boldsymbol{\mathcal{A}}^{(2)},s}$. Moreover, since $\mathrm{sc}_{s}^{(0,2)}(\mathrm{ip}^{(2,X)@}_{s}(Q))=
\mathrm{tg}_{s}^{(0,2)}(\mathrm{ip}^{(2,X)@}_{s}(P))$, i.e., these second-order paths can be $0$-composed. Thus, the interpretation of the operation symbol $\circ^{0}_{s}$ in $\mathbf{F}_{\Sigma^{\boldsymbol{\mathcal{A}}^{(2)}}}(\mathbf{Pth}_{\boldsymbol{\mathcal{A}}^{(2)}})$ is given by the interpretation of $\circ^{0}_{s}$ in $\mathbf{Pth}_{\boldsymbol{\mathcal{A}}^{(2)}}$.

Since $\mathrm{ip}^{(2,X)@}_{s}(Q\circ_{s}^{0\mathbf{T}_{\Sigma^{\boldsymbol{\mathcal{A}}^{(2)}}}(X)} P)$ is a second-order path in $\mathrm{Pth}_{\boldsymbol{\mathcal{A}}^{(2)},s}$, we have that, by Proposition~\ref{PDPT}, $Q\circ_{s}^{0\mathbf{T}_{\Sigma^{\boldsymbol{\mathcal{A}}^{(2)}}}(X)} P$ is a second-order path term in $\mathrm{PT}_{\boldsymbol{\mathcal{A}}^{(2)},s}$. That is, the operation $\circ_{s}^{0\mathbf{PT}_{\boldsymbol{\mathcal{A}}^{(2)}}}$ is well-defined.

\textsf{(vii)} For every $s\in S$ and every second-order rewrite rule $\mathfrak{p}^{(2)}\in \mathcal{A}^{(2)}_{s}$, the constant $\mathfrak{p}^{(2)\mathbf{PT}_{\boldsymbol{\mathcal{A}}^{(2)}}}$ is equal to $\mathfrak{p}^{(2)\mathbf{T}_{\Sigma^{\boldsymbol{\mathcal{A}}^{(2)}}}(X)}$, i.e., to the interpretation of $\mathfrak{p}^{(2)}$ in $\mathbf{T}_{\Sigma^{\boldsymbol{\mathcal{A}}^{(2)}}}(X)$.

Let us note that the following chain of equalities holds
\allowdisplaybreaks
\begin{align*}
\mathrm{ip}^{(2,X)@}_{s}\left(
\mathfrak{p}^{(2)\mathbf{T}_{\Sigma^{\boldsymbol{\mathcal{A}}^{(2)}}}(X)}
\right)
&=
\mathfrak{p}^{(2)\mathbf{F}_{\Sigma^{\boldsymbol{\mathcal{A}}^{(2)}}}\left(
\mathbf{Pth}_{\boldsymbol{\mathcal{A}}}^{(2)}\right)}
\tag{1}
\\&=
\mathfrak{p}^{(2)\mathbf{Pth}_{\boldsymbol{\mathcal{A}}^{(2)}}}.
\tag{2}
\end{align*}

The first equality follows from the fact that, by Definition~\ref{DDIp}, $\mathrm{ip}^{(2,X)@}$ is a $\Sigma^{\boldsymbol{\mathcal{A}}^{(2)}}$-homomorphism; finally, the second equality follows from the fact that the interpretation of the constant operation symbol $\mathfrak{p}^{(2)}$ in $\mathbf{F}_{\Sigma^{\boldsymbol{\mathcal{A}}^{(2)}}}(\mathbf{Pth}_{\boldsymbol{\mathcal{A}}^{(2)}})$ is given by the interpretation of $\mathfrak{p}^{(2)}$ in $\mathbf{Pth}_{\boldsymbol{\mathcal{A}}^{(2)}}$.

Thus, since $\mathrm{ip}^{(2,X)@}_{s}(\mathfrak{p}^{(2)\mathbf{T}_{\Sigma^{\boldsymbol{\mathcal{A}}^{(2)}}}(X)})$ is a second-order path in $\mathrm{Pth}_{\boldsymbol{\mathcal{A}}^{(2)},s}$, we have that, by  Proposition~\ref{PDPT}, $\mathfrak{p}^{(2)\mathbf{T}_{\Sigma^{\boldsymbol{\mathcal{A}}^{(2)}}}(X)}$ is a second-order path term in $\mathrm{PT}_{\boldsymbol{\mathcal{A}}^{(2)},s}$. That is, the constant $\mathfrak{p}^{(2)\mathbf{T}_{\Sigma^{\boldsymbol{\mathcal{A}}^{(2)}}}(X)}$ is well-defined.

\textsf{(viii)} For every $s\in S$, the operation $\mathrm{sc}_{s}^{1\mathbf{PT}_{\boldsymbol{\mathcal{A}}^{(2)}}}$ is equal to $\mathrm{sc}_{s}^{1\mathbf{T}_{\Sigma^{\boldsymbol{\mathcal{A}}^{(2)}}}(X)}$, i.e., to the interpretation of $\mathrm{sc}_{s}^{1}$ in $\mathbf{T}_{\Sigma^{\boldsymbol{\mathcal{A}}^{(2)}}}(X)$.

Let us note that, if $P$ is a second-order path term in $\mathrm{PT}_{\boldsymbol{\mathcal{A}}^{(2)},s}$, then the following chain of equalities holds
\allowdisplaybreaks
\begin{align*}
\mathrm{ip}^{(2,X)@}_{s}\left(
\mathrm{sc}_{s}^{1\mathbf{T}_{\Sigma^{\boldsymbol{\mathcal{A}}^{(2)}}}(X)}\left(
P
\right)\right)
&=
\mathrm{sc}_{s}^{1\mathbf{F}_{\Sigma^{\boldsymbol{\mathcal{A}}^{(2)}}}
\left(\mathbf{Pth}_{\boldsymbol{\mathcal{A}}^{(2)}}\right)}
\left(\mathrm{ip}^{(2,X)@}_{s}\left(
P
\right)\right)
\tag{1}
\\&=
\mathrm{sc}_{s}^{1\mathbf{Pth}_{\boldsymbol{\mathcal{A}}^{(2)}}}\left(
\mathrm{ip}^{(2,X)@}_{s}\left(
P
\right)\right).
\tag{2}
\end{align*}

The first equality follows from the fact that, by Definition~\ref{DDIp}, $\mathrm{ip}^{(2,X)@}$ is a $\Sigma^{\boldsymbol{\mathcal{A}}^{(2)}}$-homomorphism; finally, the second equality follows from the fact that $P$ is a second-order path term in $\mathrm{PT}_{\boldsymbol{\mathcal{A}}^{(2)},s}$. Hence, by Proposition~\ref{PDPT} $\mathrm{ip}^{(2,X)@}_{s}(P)$ is a second-order path in $\mathrm{Pth}_{\boldsymbol{\mathcal{A}}^{(2)},s}$. Thus, the interpretation of the operation symbol $\mathrm{sc}^{1}_{s}$ in $\mathbf{F}_{\Sigma^{\boldsymbol{\mathcal{A}}^{(2)}}}(\mathbf{Pth}_{\boldsymbol{\mathcal{A}}^{(2)}})$ is given by the interpretation of $\mathrm{sc}^{1}_{s}$ in $\mathbf{Pth}_{\boldsymbol{\mathcal{A}}^{(2)}}$.

Therefore, since $\mathrm{ip}^{(2,X)@}_{s}(\mathrm{sc}_{s}^{1\mathbf{T}_{\Sigma^{\boldsymbol{\mathcal{A}}^{(2)}}}(X)}(P))$ is a second-order path in $\mathrm{Pth}_{\boldsymbol{\mathcal{A}}^{(2)},s}$, we have that, by Proposition~\ref{PDPT}, $\mathrm{sc}_{s}^{1\mathbf{T}_{\Sigma^{\boldsymbol{\mathcal{A}}^{(2)}}}(X)}(P)$ is a second-order path term in $\mathrm{PT}_{\boldsymbol{\mathcal{A}}^{(2)},s}$. That is, the operation $\mathrm{sc}_{s}^{1\mathbf{PT}_{\boldsymbol{\mathcal{A}}^{(2)}}}$ is well-defined.

\textsf{(ix)} For every $s\in S$, the operation $\mathrm{tg}_{s}^{1\mathbf{PT}_{\boldsymbol{\mathcal{A}}^{(2)}}}$ is equal to $\mathrm{tg}_{s}^{1\mathbf{T}_{\Sigma^{\boldsymbol{\mathcal{A}}^{(2)}}}(X)}$, i.e., to the interpretation of $\mathrm{tg}_{s}^{1}$ in $\mathbf{T}_{\Sigma^{\boldsymbol{\mathcal{A}}^{(2)}}}(X)$.  The operation $\mathrm{tg}_{s}^{1\mathbf{PT}_{\boldsymbol{\mathcal{A}}^{(2)}}}$ is well-defined by a similar argument to that of $\mathrm{sc}_{s}^{1\mathbf{PT}_{\boldsymbol{\mathcal{A}}^{(2)}}}$.

\textsf{(x)} For every $s\in S$, the partial binary operation $\circ_{s}^{1\mathbf{PT}_{\boldsymbol{\mathcal{A}}^{(2)}}}$ is defined on second-order path terms $Q,P\in\mathrm{PT}_{\boldsymbol{\mathcal{A}}^{(2)},s}$ if, and only if, the following equation holds
$$
\mathrm{sc}_{s}^{([1],2)}\left(
\mathrm{ip}^{(2,X)@}_{s}\left(
Q
\right)\right)=
\mathrm{tg}_{s}^{([1],2)}\left(
\mathrm{ip}^{(2,X)@}_{s}\left(
P
\right)\right).
$$
For this case, $\circ_{s}^{1\mathbf{PT}_{\boldsymbol{\mathcal{A}}^{(2)}}}$ is equal to $\circ_{s}^{1\mathbf{T}_{\Sigma^{\boldsymbol{\mathcal{A}}^{(2)}}}(X)}$, i.e., to the interpretation of $\circ^{1}_{s}$ in $\mathbf{T}_{\Sigma^{\boldsymbol{\mathcal{A}}^{(2)}}}(X)$. 

Let us note that, if $Q$ and $P$ are second-order path terms in $\mathrm{PT}_{\boldsymbol{\mathcal{A}}^{(2)},s}$ with 
$$\mathrm{sc}_{s}^{([1],2)}\left(
\mathrm{ip}^{(2,X)@}_{s}\left(
Q
\right)\right)=
\mathrm{tg}_{s}^{([1],2)}\left(
\mathrm{ip}^{(2,X)@}_{s}\left(
P
\right)\right),$$ 
then the following chain of equalities holds
\begin{flushleft}
$\mathrm{ip}^{(2,X)@}_{s}\left(
Q\circ_{s}^{1\mathbf{T}_{\Sigma^{\boldsymbol{\mathcal{A}}^{(2)}}}(X)} P
\right)$
\allowdisplaybreaks
\begin{align*}
\qquad
&=
\left(\mathrm{ip}^{(2,X)@}_{s}\left(
Q
\right)\right)
\circ_{s}^{1\mathbf{F}_{\Sigma^{\boldsymbol{\mathcal{A}}^{(2)}}}
\left(\mathbf{Pth}_{\boldsymbol{\mathcal{A}}^{(2)}}\right)}
\left(\mathrm{ip}^{(2,X)@}_{s}\left(
P
\right)\right)
\tag{1}
\\&=
\left(\mathrm{ip}^{(2,X)@}_{s}\left(Q
\right)\right)
\circ_{s}^{1\mathbf{Pth}_{\boldsymbol{\mathcal{A}}^{(2)}}}
\left(\mathrm{ip}^{(2,X)@}_{s}\left(
P
\right)\right).
\tag{2}
\end{align*}
\end{flushleft}

The first equality follows from the fact that, by Definition~\ref{DDIp}, $\mathrm{ip}^{(2,X)@}$ is a $\Sigma^{\boldsymbol{\mathcal{A}}^{(2)}}$-homomorphism; finally, the second equality follows from the fact that $Q$ and $P$ are second-order path terms in $\mathrm{PT}_{\boldsymbol{\mathcal{A}}^{(2)},s}$. Hence, in virtue of Proposition~\ref{PDPT}  $\mathrm{ip}^{(2,X)@}_{s}(Q)$ and $\mathrm{ip}^{(2,X)@}_{s}(P)$ are second-order paths in $\mathrm{Pth}_{\boldsymbol{\mathcal{A}}^{(2)},s}$. Moreover, $\mathrm{sc}_{s}^{([1],2)}(\mathrm{ip}^{(2,X)@}_{s}(Q))=
\mathrm{tg}_{s}^{([1],2)}(\mathrm{ip}^{(2,X)@}_{s}(P))$, i.e., these second-order paths can be $1$-composed. Thus, the interpretation of the operation symbol $\circ^{1}_{s}$ in $\mathbf{F}_{\Sigma^{\boldsymbol{\mathcal{A}}^{(2)}}}(\mathbf{Pth}_{\boldsymbol{\mathcal{A}}^{(2)}})$ is given by the interpretation of $\circ^{1}_{s}$ in $\mathbf{Pth}_{\boldsymbol{\mathcal{A}}^{(2)}}$.

Since $\mathrm{ip}^{(2,X)@}_{s}(Q\circ_{s}^{1\mathbf{T}_{\Sigma^{\boldsymbol{\mathcal{A}}^{(2)}}}(X)} P)$ is a second-order path in $\mathrm{Pth}_{\boldsymbol{\mathcal{A}}^{(2)},s}$, we have that, by Proposition~\ref{PDPT}, $Q\circ_{s}^{1\mathbf{T}_{\Sigma^{\boldsymbol{\mathcal{A}}^{(2)}}}(X)} P$ is a second-order path term in $\mathrm{PT}_{\boldsymbol{\mathcal{A}}^{(2)},s}$. That is, the operation $\circ_{s}^{1\mathbf{PT}_{\boldsymbol{\mathcal{A}}^{(2)}}}$ is well-defined.

This completes the definition of the partial $\Sigma^{\boldsymbol{\mathcal{A}}}$-algebra $\mathbf{PT}_{\boldsymbol{\mathcal{A}}^{(2)}}$.
\end{proof}

We next prove that the many-sorted partial $\Sigma^{\boldsymbol{\mathcal{A}}^{(2)}}$-algebra $\mathbf{PT}_{\boldsymbol{\mathcal{A}}^{(2)}}$ is a partial Dedekind-Peano $\Sigma^{\boldsymbol{\mathcal{A}}^{(2)}}$-algebra. 

\begin{restatable}{proposition}{PDPTPDP}
\label{PDPTPDP} 
The many-sorted partial $\Sigma^{\boldsymbol{\mathcal{A}}^{(2)}}$-algebra $\mathbf{PT}_{\boldsymbol{\mathcal{A}}^{(2)}}$ is a partial Dedekind-Peano $\Sigma^{\boldsymbol{\mathcal{A}}^{(2)}}$-algebra.
\end{restatable}
\begin{proof} 
According to Definition~\ref{DPAlgDP}, we need to check that the following items hold

\textsf{PDP1.} For every $(\mathbf{s},s)\in S^{\star}\times S$ and every $\gamma\in\Sigma^{\boldsymbol{\mathcal{A}}^{(2)}}_{\mathbf{s},s}$, the operation 
$$\gamma^{\mathbf{PT}_{\boldsymbol{\mathcal{A}}^{(2)}}}\colon
\mathbf{PT}_{\boldsymbol{\mathcal{A}}^{(2)},\mathbf{s}}
\dmor
\mathbf{PT}_{\boldsymbol{\mathcal{A}}^{(2)},s}
$$
is injective, i.e., for every $(Q_{j})_{j\in\bb{\mathbf{s}}}, (P_{j})_{j\in\bb{\mathbf{s}}}$ in $\mathrm{Dom}(\gamma)$, if 
$$\gamma^{\mathbf{PT}_{\boldsymbol{\mathcal{A}}^{(2)}}}
\left(\left(
Q_{j}
\right)_{j\in\bb{\mathbf{s}}}\right)
=\gamma^{\mathbf{PT}_{\boldsymbol{\mathcal{A}}^{(2)}}}
\left(\left(
P_{j}
\right)_{j\in\bb{\mathbf{s}}}\right),$$
then $(Q_{j})_{j\in\bb{\mathbf{s}}}=(P_{j})_{j\in\bb{\mathbf{s}}}$.

\textsf{PDP2.} For every $s\in S$ and every $\tau,\gamma\in\Sigma^{\boldsymbol{\mathcal{A}}^{(2)}}_{\cdot, s}$, if $\tau\neq \gamma$, then 
$$
\tau^{\mathbf{PT}_{\boldsymbol{\mathcal{A}}^{(2)}}}\left[
\mathrm{Dom}\left(
\tau^{\mathbf{PT}_{\boldsymbol{\mathcal{A}}^{(2)}}}
\right)\right]
\cap
\gamma^{\mathbf{PT}_{\boldsymbol{\mathcal{A}}^{(2)}}}\left[
\mathrm{Dom}\left(
\gamma^{\mathbf{PT}_{\boldsymbol{\mathcal{A}}^{(2)}}}
\right)\right]
=
\varnothing.
$$

\textsf{PDP3.} The following equality holds
\[
\textstyle
\mathrm{Sg}_{\mathbf{PT}_{\boldsymbol{\mathcal{A}}^{(2)}}}
\left(
\mathrm{PT}_{\boldsymbol{\mathcal{A}}^{(2)}}
-
\left(
\bigcup_{\gamma\in\Sigma^{\boldsymbol{\mathcal{A}}^{(2)}}_{\cdot,s}}
\mathrm{Im}\left(
\gamma^{\mathbf{PT}_{\boldsymbol{\mathcal{A}}^{(2)}}}
\right)
\right)_{s\in S}
\right)
=
\mathrm{PT}_{\boldsymbol{\mathcal{A}}^{(2)}}.
\]

This is true, because this is the case in the many-sorted free $\Sigma^{\boldsymbol{\mathcal{A}}^{(2)}}$-algebra $\mathbf{T}_{\Sigma^{\boldsymbol{\mathcal{A}}^{(2)}}}(X)$. Note that the basis of Dedekind-Peano of $\mathbf{PT}_{\boldsymbol{\mathcal{A}}^{(2)}}$ is given by
$$
\textstyle
\mathrm{PT}_{\boldsymbol{\mathcal{A}}^{(2)}}
-
\left(
\bigcup_{\gamma\in\Sigma^{\boldsymbol{\mathcal{A}}^{(2)}}_{\cdot,s}}
\mathrm{Im}\left(
\gamma^{\mathbf{PT}_{\boldsymbol{\mathcal{A}}^{(2)}}}
\right)
\right)_{s\in S}=X.
$$
That is, $X$ is the basis of Dedekind-Peano of $\mathbf{PT}_{\boldsymbol{\mathcal{A}}^{(2)}}$.
\end{proof}

\section{
\texorpdfstring
{An Artinian partial order on $\coprod\mathrm{PT}_{\boldsymbol{\mathcal{A}}^{(2)}}$}
{An Artinian order on second-order path terms}
}

In this section we consider a partial order on $\coprod\mathrm{PT}_{\boldsymbol{\mathcal{A}}^{(2)}}$, the coproduct of the many-sorted set of second-order path terms.

The next corollary states that the subterms of a second-order path term are also second-order path terms.

\begin{restatable}{corollary}{CDPTSubt}
\label{CDPTSubt} Let $s$ be a sort in $S$ and $P$ a second-order path term in $\mathrm{PT}_{\boldsymbol{\mathcal{A}}^{(2)},s}$, then $\mathrm{Subt}_{\Sigma^{\boldsymbol{\mathcal{A}}^{(2)}}}(P)\subseteq \mathrm{PT}_{\boldsymbol{\mathcal{A}}^{(2)}}.$
\end{restatable}
\begin{proof}
It follows from Definition~\ref{DDPT} and Lemma~\ref{LDThetaCongSub}.
\end{proof}

As a consequence, we can consider the restriction of the Artinian partial order $\leq_{\mathbf{T}_{\Sigma^{\boldsymbol{\mathcal{A}}^{(2)}}}(X)}$ on $\coprod\mathrm{T}_{\Sigma^{\boldsymbol{\mathcal{A}}^{(2)}}}$ to the coproduct of the many-sorted set of second-order path terms.

\begin{restatable}{definition}{DDPTOrd}
\label{DDPTOrd}
\index{partial order!second-order!$\leq_{\mathbf{PT}_{\boldsymbol{\mathcal{A}}^{(2)}}}$}
Let $\leq_{\mathbf{PT}_{\boldsymbol{\mathcal{A}}^{(2)}}}$ be the binary relation on $\mathrm{PT}_{\boldsymbol{\mathcal{A}}^{(2)}}$ containing every pair $((Q,t),(P,s))$ in $(\coprod\mathrm{PT}_{\boldsymbol{\mathcal{A}}^{(2)}})^{2}$ satisfying that
$$
\left(
Q,t
\right)
\leq_{\mathbf{T}_{\Sigma^{\boldsymbol{\mathcal{A}}^{(2)}}}(X)}
\left(
P,s
\right).
$$

That is, $(Q,t)$ $\leq_{\mathbf{PT}_{\boldsymbol{\mathcal{A}}^{(2)}}}$-precedes $(P,s)$ if, and only if, $Q$ is a subterm of type $t$ of $P$.
\end{restatable}

\begin{restatable}{proposition}{PDPTOrdArt}
\label{PDPTOrdArt}
$(\coprod\mathrm{PT}_{\boldsymbol{\mathcal{A}}^{(2)}}, \leq_{\mathbf{PT}_{\boldsymbol{\mathcal{A}}^{(2)}}})$ is a partially ordered set. Moreover, in this partially ordered set there is not any strictly decreasing $\omega_{0}$-chain, i.e., $(\coprod\mathrm{PT}_{\boldsymbol{\mathcal{A}}^{(2)}}, \leq_{\mathbf{PT}_{\boldsymbol{\mathcal{A}}^{(2)}}})$ is an Artinian partially ordered set.
\end{restatable}

\section{
\texorpdfstring
{On the quotient $\llbracket\mathrm{PT}_{\boldsymbol{\mathcal{A}}^{(2)}}\rrbracket$}
{On the quotient of second-order path terms}
}

In this section we introduce the $S$-sorted set of second-order path term classes. This many-sorted set is defined as the quotient of the $S$-sorted set of second-order path terms by the corresponding restriction of $\Theta^{\llbracket 2 \rrbracket}$,
 the $\Sigma^{\boldsymbol{\mathcal{A}}^{(2)}}$-congruence on $\mathbf{T}_{\Sigma^{\boldsymbol{\mathcal{A}}^{(2)}}}(X)$ introduced in Definition~\ref{DDVCong}.
 
\begin{restatable}{convention}{CDPTClass}
\label{CDPTClass} 
\index{path-term!second-order!$\llbracket P\rrbracket_{s}$}
\index{path-term!second-order!$\llbracket\mathrm{PT}_{\boldsymbol{\mathcal{A}}^{(2)}}\rrbracket$}
To simplify the notation, for a sort $s\in S$ and a second-order path term $P\in\mathrm{PT}_{\boldsymbol{\mathcal{A}}^{(2)},s}$, we will let $\llbracket P\rrbracket_{s}$ stand for $\llbracket P\rrbracket_{\Theta^{\llbracket 2 \rrbracket}_{s}}$, the 
$\Theta^{\llbracket 2\rrbracket}_{s}$-equivalence class of $P$, and we will call it the \emph{second-order path term class of} $P$. 
\end{restatable}

\begin{restatable}{definition}{DDPTQuot}
\label{DDPTQuot} 
\index{path terms!second-order!$\llbracket\mathrm{PT}_{\boldsymbol{\mathcal{A}}^{(2)}}\rrbracket$}
We denote by $\llbracket\mathrm{PT}_{\boldsymbol{\mathcal{A}}^{(2)}}\rrbracket$ the image of $\mathrm{PT}_{\boldsymbol{\mathcal{A}}^{(2)}}$ under the projection for the $\Sigma^{\boldsymbol{\mathcal{A}}^{(2)}}$-congruence $\Theta^{\llbracket 2 \rrbracket}$. We call it the $S$-sorted set of \emph{second-order path term classes}. Note that $\llbracket\mathrm{PT}_{\boldsymbol{\mathcal{A}}^{(2)}}\rrbracket$  is a subset of the quotient $\mathrm{T}_{\Sigma^{\boldsymbol{\mathcal{A}}^{(2)}}}(X)/{\Theta^{\llbracket 2 \rrbracket}}$. Actually, we have that
$$
\llbracket\mathrm{PT}_{\boldsymbol{\mathcal{A}}^{(2)}}\rrbracket = 
\mathrm{PT}_{\boldsymbol{\mathcal{A}}^{(2)}}/\Theta^{\llbracket 2 \rrbracket}\!\upharpoonright\!\mathrm{PT}_{\boldsymbol{\mathcal{A}}^{(2)}}.
$$

\index{projection!second-order!$\mathrm{pr}^{\Theta^{\llbracket 2 \rrbracket}}$}
As expected, the projection $\mathrm{pr}^{\Theta^{\llbracket 2 \rrbracket}}$ from $\mathrm{T}_{\Sigma^{\boldsymbol{\mathcal{A}}^{(2)}}}(X)$ to $\mathrm{T}_{\Sigma^{\boldsymbol{\mathcal{A}}^{(2)}}}(X)/{\Theta^{\llbracket 2 \rrbracket}}$ birrestricts to $\mathrm{PT}_{\boldsymbol{\mathcal{A}}^{(2)}}$ and
$\llbracket\mathrm{PT}_{\boldsymbol{\mathcal{A}}^{(2)}}\rrbracket$, respectively. To avoid further notation, we resignify the already existing notation. 

We let $\mathrm{pr}^{\Theta^{\llbracket 2 \rrbracket}}$ be the projection from $\mathrm{PT}_{\boldsymbol{\mathcal{A}}^{(2)}}$ to $\llbracket\mathrm{PT}_{\boldsymbol{\mathcal{A}}^{(2)}}\rrbracket$ that, for every sort $s\in S$, maps a second-order path term $P$ in $\mathrm{PT}_{\boldsymbol{\mathcal{A}}^{(2)},s}$ to $\llbracket P \rrbracket_{s}$, its equivalence class under $\Theta^{\llbracket 2 \rrbracket}_{s}$.
$$
\mathrm{pr}^{\Theta^{\llbracket 2 \rrbracket}}
\colon
\mathrm{PT}_{\boldsymbol{\mathcal{A}}^{(2)}}
\mor
\llbracket\mathrm{PT}_{\boldsymbol{\mathcal{A}}^{(2)}}\rrbracket.
$$
\end{restatable}

We next introduce the notation we will later use for the composition of the embeddings introduced in Definition~\ref{DDEta} with the projection to $\Theta^{\llbracket 2 \rrbracket}$.

\begin{restatable}{definition}{DDPTQEta}
\label{DDPTQEta}  Let $X$ be an $S$-sorted set and let $\llbracket\mathrm{PT}_{\boldsymbol{\mathcal{A}}^{(2)}}\rrbracket$ be the many-sorted set of second-order path term classes.  We will denote by
\begin{enumerate}
\item $\eta^{(\llbracket 2\rrbracket,X)}$ the $S$-sorted mapping from $X$ to 
$\llbracket\mathrm{PT}_{\boldsymbol{\mathcal{A}}^{(2)}}\rrbracket$ given by the composition $\eta^{(\llbracket 2\rrbracket,X)}=\mathrm{pr}^{\Theta^{\llbracket 2 \rrbracket}}\circ\eta^{(2,X)}$, i.e., for every sort $s\in S$,  $\eta^{(\llbracket 2\rrbracket,X)}$ sends a variable $x\in X_{s}$ to the class $\llbracket \eta^{(2,X)}_{s}(x)\rrbracket_{s}$ in $\llbracket\mathrm{PT}_{\boldsymbol{\mathcal{A}}^{(2)}}\rrbracket_{s}$.
\index{inclusion!second-order!$\eta^{(\llbracket 2\rrbracket,X)}$}
\item $\eta^{(\llbracket 2\rrbracket,\mathcal{A})}$ the $S$-sorted mapping from $\mathcal{A}$ to 
$\llbracket\mathrm{PT}_{\boldsymbol{\mathcal{A}}^{(2)}}\rrbracket$ given by the composition $\eta^{(\llbracket 2\rrbracket,\mathcal{A})}=\mathrm{pr}^{\Theta^{\llbracket 2 \rrbracket}}\circ\eta^{(2,\mathcal{A})}$, i.e., for every sort $s\in S$,  $\eta^{(\llbracket 2\rrbracket,\mathcal{A})}$ sends a rewrite rule  $\mathfrak{p}\in \mathcal{A}_{s}$ to the class $\llbracket \eta^{(2,\mathcal{A})}_{s}(\mathfrak{p})\rrbracket_{s}$ in $\llbracket\mathrm{PT}_{\boldsymbol{\mathcal{A}}^{(2)}}\rrbracket_{s}$.
\index{inclusion!second-order!$\eta^{(\llbracket 2\rrbracket,\mathcal{A})}$}
\item $\eta^{(\llbracket 2\rrbracket,\mathcal{A}^{(2)})}$ the $S$-sorted mapping from $\mathcal{A}^{(2)}$ to 
$\llbracket\mathrm{PT}_{\boldsymbol{\mathcal{A}}^{(2)}}\rrbracket$ given by the composition $\eta^{(\llbracket 2\rrbracket,\mathcal{A}^{(2)})}=\mathrm{pr}^{\Theta^{\llbracket 2 \rrbracket}}\circ\eta^{(2,\mathcal{A}^{(2)})}$, i.e., for every sort $s\in S$,  $\eta^{(\llbracket 2\rrbracket,\mathcal{A}^{(2)})}$ sends a second-order rewrite rule  $\mathfrak{p}^{(2)}\in \mathcal{A}^{(2)}_{s}$ to the class $\llbracket\eta^{(2,\mathcal{A}^{(2)})}_{s}(\mathfrak{p}^{(2)})\rrbracket_{s}$ in $\llbracket\mathrm{PT}_{\boldsymbol{\mathcal{A}}^{(2)}}\rrbracket_{s}$.
\index{inclusion!second-order!$\eta^{(\llbracket 2\rrbracket,\mathcal{A}^{(2)})}$}
\end{enumerate}
The above $S$-sorted mappings are depicted in the diagram of Figure~\ref{FDPTQEta}.
\end{restatable}

\begin{figure}
\begin{tikzpicture}
[ACliment/.style={-{To [angle'=45, length=5.75pt, width=4pt, round]}},scale=.8]
\node[] (x) at (0,0) [] {$X$};
\node[] (a) at (0,-1.5) [] {$\mathcal{A}$};
\node[] (a2) at (0,-3) [] {$\mathcal{A}^{(2)}$};
\node[] (T) at (6,-3) [] {$
\llbracket\mathrm{PT}_{\boldsymbol{\mathcal{A}}^{(2)}}\rrbracket$};
\draw[ACliment, bend left=20]  (x) to node [above right] {$\eta^{(\llbracket 2\rrbracket,X)}$} (T);
\draw[ACliment, bend left=10]  (a) to node [midway, fill=white] {$\eta^{(\llbracket 2\rrbracket,\mathcal{A})}$} (T);
\draw[ACliment]  (a2) to node [below] {$\eta^{(\llbracket 2\rrbracket,\mathcal{A}^{(2)})}$} (T);
\end{tikzpicture}
\caption{Quotient term embeddings relative to $X$, $\mathcal{A}$ and $\mathcal{A}^{(2)}$ at layer 2.}\label{FDPTQEta}
\end{figure}

We now consider the embeddings introduced in Propositions~\ref{PDUEmb} and~\ref{PDZEmb} to second-order path term classes. 

\begin{restatable}{definition}{DDPTQEtas}
\label{DDPTQEtas}  Let $X$ be an $S$-sorted set and let $\llbracket\mathrm{PT}_{\boldsymbol{\mathcal{A}}^{(2)}}\rrbracket$ be the many-sorted set of second-order path term classes.  We will denote by
\begin{enumerate}
\item  $\eta^{(\llbracket 2\rrbracket,[1])\sharp}$ the $S$-sorted mapping from $[\mathrm{PT}_{\boldsymbol{\mathcal{A}}}]$ to 
$\llbracket\mathrm{PT}_{\boldsymbol{\mathcal{A}}^{(2)}}\rrbracket$ such that, for every sort $s\in S$, sends a path term class $[P]_{s}$ in $[\mathrm{PT}_{\boldsymbol{\mathcal{A}}}]_{s}$ to the second-order path term class  $
\llbracket \eta^{(2,1)\sharp}_{s}(P)\rrbracket_{s}$ in $\llbracket\mathrm{PT}_{\boldsymbol{\mathcal{A}}^{(2)}}\rrbracket_{s}$. Note that this $S$-sorted mapping is well-defined since, according to Proposition~\ref{PDThetaCongDU}, $\eta^{(2,1)\sharp^{2}}[\Theta^{[1]}]\subseteq\Theta^{[2]}\subseteq\Theta^{\llbracket 2 \rrbracket}$.
\index{inclusion!second-order!$\eta^{(\llbracket 2\rrbracket,[1])\sharp}$}
\item $\eta^{(\llbracket 2\rrbracket,0)\sharp}$ the $S$-sorted mapping from $\mathrm{T}_{\Sigma}(X)$ to 
$\llbracket\mathrm{PT}_{\boldsymbol{\mathcal{A}}^{(2)}}\rrbracket$ such that, for every sort $s\in S$, sends a term $P\in \mathrm{T}_{\Sigma}(X)_{s}$ to the second-order path term class  $\llbracket \eta^{(2,0)\sharp}_{s}(P)\rrbracket_{s}$ in $\llbracket\mathrm{PT}_{\boldsymbol{\mathcal{A}}^{(2)}}\rrbracket_{s}$.
\index{inclusion!second-order!$\eta^{(\llbracket 2\rrbracket,0)\sharp}$}
\end{enumerate}
\end{restatable}

\begin{proposition}\label{PDPTQX} The following equalities hold
\begin{multicols}{2}
\begin{itemize}
\item[(i)] $\eta^{(\llbracket 2\rrbracket,X)}=\eta^{(\llbracket 2\rrbracket,[1])\sharp}\circ\eta^{([1],X)}$;
\item[(ii)] $\eta^{(\llbracket 2\rrbracket,0)\sharp}=\eta^{(\llbracket 2\rrbracket,[1])\sharp}\circ\eta^{([1],0)\sharp}.$ 
\end{itemize}
\end{multicols}
\end{proposition}
\begin{proof}
This proposition entails that the diagram in Figure~\ref{FDPTQXEq} commutes.

Let $s$ be a sort in $S$ and $x$ a variable in $X_{s}$. 
The following chain of equalities holds
\allowdisplaybreaks
\begin{align*}
\eta^{(\llbracket 2\rrbracket,[1])\sharp}_{s}\left(
\eta^{([1],X)}_{s}\left(
x
\right)
\right)
&=
\eta^{(\llbracket 2\rrbracket,[1])\sharp}_{s}\left(
\left[
\eta^{(1,X)}_{s}(x)
\right]_{s}
\right)
\tag{1}
\\&=
\biggl\llbracket
\eta^{(2,1)\sharp}_{s}\left(
\eta^{(1,X)}_{s}\left(
x
\right)
\right)
\biggr\rrbracket_{s}
\tag{2}
\\&=
\Bigl\llbracket
\eta^{(2,X)}_{s}\left(
x
\right)
\Bigr\rrbracket_{s}
\tag{3}
\\&=
\eta^{(\llbracket 2\rrbracket,X)}_{s}\left(
x
\right).
\tag{4}
\end{align*}

The first equality applies the mapping $\eta^{([1],X)}$ according to Definition~\ref{DPTQEta}; the second equality applies the mapping $\eta^{(\llbracket 2\rrbracket,[1])\sharp}$ according to Definition~\ref{DDPTQEtas}; the third equality follows from Proposition~\ref{PDUEmb}; finally, the last equality recovers the mapping $\eta^{(\llbracket 2\rrbracket,X)}$ introduced in Proposition~\ref{DDPTQEta}.

Regarding the second equality, let $s$ be a sort in $S$ and $P$ a term  in $\mathrm{T}_{\Sigma}(X)_{s}$. 
The following chain of equalities holds
\allowdisplaybreaks
\begin{align*}
\eta^{(\llbracket 2\rrbracket,[1])\sharp}_{s}\left(
\eta^{([1],0)\sharp}_{s}\left(
P
\right)
\right)
&=
\eta^{(\llbracket 2\rrbracket,[1])\sharp}_{s}\left(
\left[
\eta^{(1,0)\sharp}_{s}\left(
P
\right)
\right]_{s}
\right)
\tag{1}
\\&=
\biggl\llbracket
\eta^{(2,1)\sharp}_{s}\left(
\eta^{(1,0)\sharp}_{s}\left(
P
\right)
\right)
\biggr\rrbracket_{s}
\tag{2}
\\&=
\Bigl\llbracket
\eta^{(2,0)\sharp}_{s}\left(
P
\right)
\Bigr\rrbracket_{s}
\tag{3}
\\&=
\eta^{(\llbracket 2\rrbracket,0)\sharp}_{s}\left(
P
\right).
\tag{4}
\end{align*}

The first equality applies the mapping $\eta^{([1],0)\sharp}$ according to Definition~\ref{DPTQEta}; the second equality applies the mapping $\eta^{(\llbracket 2\rrbracket,[1])\sharp}$ according to Definition~\ref{DDPTQEtas}; the third equality follows from Proposition~\ref{PDEmb}; finally, the last equaliy recovers the mapping $\eta^{(\llbracket 2\rrbracket,0)\sharp}$ according to Definition~\ref{DDPTQEtas}.

This concludes the proof.
\end{proof}

\begin{figure}
\begin{center}
\begin{tikzpicture}
[ACliment/.style={-{To [angle'=45, length=5.75pt, width=4pt, round]}},scale=0.8]
\node[] (x) at (0,0) [] {$X$};
\node[] (t) at (6,0) [] {$\mathrm{T}_{\Sigma}(X)$};
\node[] (t1) at (6,-2) [] {$[\mathrm{PT}_{\boldsymbol{\mathcal{A}}}]$};
\node[] (t2) at (6,-4) [] {$\llbracket\mathrm{PT}_{\boldsymbol{\mathcal{A}}^{(2)}}\rrbracket$};

\draw[ACliment]  (x) to node [ above right]
{$\textstyle \eta^{(0,X)}$} (t);
\draw[ACliment, bend right=10]  (x) to node [midway, fill=white]
{$\textstyle\eta^{([1],X)}$} (t1);
\draw[ACliment, bend right=20]  (x) to node [ below left]
{$\textstyle \eta^{(\llbracket 2\rrbracket,X)}$} (t2);

\draw[ACliment] 
(t) to node [right] {$\textstyle \eta^{([1],0)\sharp}$}  (t1);
\draw[ACliment] 
(t1) to node [right] {$\textstyle \eta^{(\llbracket 2\rrbracket,[1])\sharp}$}  (t2);

\draw[ACliment, rounded corners] (t.east)
--
 ($(t.east)+(2,0)$)
-- node [right] {$\textstyle \eta^{(\llbracket 2\rrbracket,0)\sharp}$}
($(t.east)+(2,-4)$)
-- (t2.east);
\end{tikzpicture}
\end{center}
\caption{Quotient embeddings relative to $X$ at layers 0, 1 \& 2.}
\label{FDPTQX}
\end{figure}

\begin{proposition}\label{PDPTQA} The equality
$\eta^{(\llbracket 2\rrbracket,[1])\sharp}\circ\eta^{([1],\mathcal{A})}=\eta^{(\llbracket 2\rrbracket,\mathcal{A})}$
 holds.
\end{proposition}
\begin{proof}
This proposition entails that the diagram in Figure~\ref{FDPTQA} commutes. 

Let $s$ be a sort in $S$ and $\mathfrak{p}$ a rewrite rule in $\mathcal{A}_{s}$. The following chain of equalities holds
\allowdisplaybreaks
\begin{align*}
\eta^{(\llbracket 2\rrbracket,[1])\sharp}_{s}\left(
\eta^{([1],\mathcal{A})}_{s}\left(
\mathfrak{p}
\right)
\right)
&=
\eta^{(\llbracket 2\rrbracket,[1])\sharp}_{s}\left(
\left[
\mathfrak{p}
\right]_{s}
\right)
\tag{1}
\\&=
\biggl\llbracket
\eta^{(2,1)\sharp}_{s}\left(
\eta^{(1,\mathcal{A})}_{s}\left(
\mathfrak{p}
\right)
\right)
\biggr\rrbracket_{s}
\tag{(2)}
\\&=
\Bigl\llbracket
\eta^{(2,\mathcal{A})}_{s}\left(
\mathfrak{p}
\right)
\Bigr\rrbracket_{s}
\tag{3}
\\&=
\eta^{(\llbracket 2\rrbracket,\mathcal{A})}_{s}\left(
\mathfrak{p}
\right).
\tag{4}
\end{align*}

The first equality applies the mapping $\eta^{([1],\mathcal{A})}$ according to Definition~\ref{DPTQEta}; the second equality applies the mapping $\eta^{(\llbracket 2\rrbracket,[1])\sharp}$ according to Definition~\ref{DDPTQEtas}; the third equality follows from Proposition~\ref{PDAEmb}; finally, the last equality recovers the mapping $\eta^{(\llbracket 2\rrbracket,\mathcal{A})}$ introduced in Definition~\ref{DDPTQEta}.

This completes the proof.
\end{proof}

\begin{figure}
\begin{center}
\begin{tikzpicture}
[ACliment/.style={-{To [angle'=45, length=5.75pt, width=4pt, round]}},scale=0.8]
\node[] (x) at (0,0) [] {$\mathcal{A}$};
\node[] (t1) at (6,0) [] {$[\mathrm{PT}_{\boldsymbol{\mathcal{A}}}]$};
\node[] (t2) at (6,-2) [] {$\llbracket\mathrm{PT}_{\boldsymbol{\mathcal{A}}^{(2)}}\rrbracket$};

\draw[ACliment]  (x) to node [above]
{$\textstyle \eta^{([1],\mathcal{A})}$} (t1);
\draw[ACliment, bend right=10]  (x) to node [ below left]
{$\textstyle \eta^{(\llbracket 2\rrbracket,\mathcal{A})}$} (t2);

\draw[ACliment] 
(t1) to node [right] {$\textstyle \eta^{(\llbracket 2\rrbracket,[1])\sharp}$}  (t2);

\end{tikzpicture}
\end{center}
\caption{Quotient embeddings relative to $\mathcal{A}$ at layers 1 \& 2.}
\label{FDPTQA}
\end{figure}

We next consider the case of the  $\mathrm{ip}^{(2,X)@}$ mapping.

\begin{restatable}{corollary}{CDPTQKer}
\label{CDPTQKer} 
Let $s$ be a sort in $S$, and $Q$, $P$ second-order path terms in $\mathrm{PT}_{\boldsymbol{\mathcal{A}}^{(2)},s}$ such that $(Q,P)\in\Theta^{\llbracket 2 \rrbracket}_{s}$. Then 
$
(Q,P)\in\mathrm{Ker}(\mathrm{pr}^{\llbracket \cdot \rrbracket}\circ\mathrm{ip}^{(2,X)@}).
$
\end{restatable}
\begin{proof}
It follows from Corollary~\ref{CDVCong} and Lemma~\ref{LDVCong}.
\end{proof}

\begin{restatable}{definition}{DDPTQIp}
\label{DDPTQIp} 
\index{identity!second-order!$\mathrm{ip}^{(\llbracket 2\rrbracket,X)a}$}
We let $\mathrm{ip}^{(\llbracket 2\rrbracket,X)@}$ stand for 
$$\left(\mathrm{pr}^{\llbracket \cdot \rrbracket}\circ\mathrm{ip}^{(2,X)@}\right)^{\mathrm{m}}\circ\mathrm{p}^{
\mathrm{pr}^{\llbracket \cdot \rrbracket}\circ\mathrm{ip}^{(2,X)@},
\Theta^{\llbracket 2 \rrbracket}
},$$ 
the unique many-sorted mapping from $\llbracket\mathrm{PT}_{\boldsymbol{\mathcal{A}}^{(2)}}\rrbracket$ to $\llbracket\mathrm{Pth}_{\boldsymbol{\mathcal{A}}^{(2)}}\rrbracket$ satisfying that 
$$
\mathrm{pr}^{\llbracket \cdot \rrbracket}\circ\mathrm{ip}^{(2,X)@}
=
\mathrm{p}^{
\mathrm{pr}^{\llbracket \cdot \rrbracket}\circ\mathrm{ip}^{(2,X)@},
\Theta^{\llbracket 2 \rrbracket}
}
\circ
\mathrm{pr}^{\Theta^{\llbracket 2 \rrbracket}}.
$$
This $S$-sorted mapping is well-defined according to Corollary~\ref{CDPTQKer}.
\end{restatable}

We conclude with the second-order Curry-Howard mapping.

\begin{restatable}{definition}{DDPTQDCH}
\label{DDPTQDCH} 
\index{Curry-Howard!second-order!$\mathrm{CH}^{\llbracket 2\rrbracket}$}
We let $\mathrm{CH}^{\llbracket 2\rrbracket}$ stand for the $S$-sorted mapping from $\llbracket \mathrm{Pth}_{\boldsymbol{\mathcal{A}}^{(2)}}\rrbracket $ to $\llbracket\mathrm{PT}_{\boldsymbol{\mathcal{A}}^{(2)}}\rrbracket$ given by the composition $ \mathrm{pr}^{\Theta^{\llbracket 2 \rrbracket}}\circ\mathrm{CH}^{(2)\mathrm{m}}$, i.e.,
$$
\mathrm{CH}^{\llbracket 2\rrbracket}\colon \llbracket \mathrm{Pth}_{\boldsymbol{\mathcal{A}}^{(2)}}\rrbracket\mor
\llbracket\mathrm{PT}_{\boldsymbol{\mathcal{A}}^{(2)}}\rrbracket.
$$

That this mapping is well-defined follows from Proposition~\ref{PDVDCH}.
\end{restatable}

\begin{proposition}\label{PDPTQXEq} The following equalities holds
\begin{multicols}{2}
\begin{itemize}
\item[(i)] $\mathrm{ip}^{(\llbracket 2\rrbracket,X)@}\circ\eta^{(\llbracket 2\rrbracket,X)}=\mathrm{ip}^{(\llbracket 2\rrbracket,X)};$
\item[(ii)] $\mathrm{ip}^{(\llbracket 2\rrbracket,X)@}\circ\eta^{(\llbracket 2\rrbracket,0)\sharp}=\mathrm{ip}^{(\llbracket 2\rrbracket,0)\sharp};$
\item[(iii)] $\mathrm{CH}^{\llbracket 2\rrbracket}\circ\mathrm{ip}^{(\llbracket 2\rrbracket,X)}=\eta^{(\llbracket 2\rrbracket,X)};$
\item[(iv)] $\mathrm{CH}^{\llbracket 2\rrbracket}\circ\mathrm{ip}^{(\llbracket 2\rrbracket,0)\sharp}=\eta^{(\llbracket 2\rrbracket,0)\sharp};$
\end{itemize}
\end{multicols}
The reader is advised to consult the diagram presented in Figure~\ref{FDPTQX}.
\end{proposition}
\begin{proof}
Regarding the first item, let $s$ be a sort in $S$ and $x$ a variable in $X_{s}$. The following chain of equalities holds
\allowdisplaybreaks
\begin{align*}
\mathrm{ip}^{(\llbracket 2\rrbracket,X)@}_{s}\left(
\eta^{(\llbracket 2\rrbracket,X)}_{s}\left(
x
\right)\right)
&=
\mathrm{ip}^{(\llbracket 2\rrbracket,X)@}_{s}\left(
\bigl\llbracket
\eta^{(2,X)}_{s}\left(
x
\right)
\bigr\rrbracket_{s}
\right)
\tag{1}
\\&=
\Bigl\llbracket
\mathrm{ip}^{(2,X)@}_{s}\left(
\eta^{(2,X)}_{s}\left(
x
\right)
\right)
\Bigr\rrbracket_{s}
\tag{2}
\\&=
\bigl\llbracket
\mathrm{ip}^{(2,X)}_{s}\left(
x
\right)
\bigr\rrbracket_{s}
\tag{3}
\\&=
\mathrm{ip}^{(\llbracket 2\rrbracket,X)}_{s}\left(
x
\right).
\tag{4}
\end{align*}

The first equality applies the mapping $\eta^{(\llbracket 2\rrbracket,X)}$ introduced in Definition~\ref{DDPTQEta}; the second equality applies the mapping $\mathrm{ip}^{(\llbracket 2\rrbracket,X)@}$ introduced in Definition~\ref{DDPTQIp}; the third equality follows from Definition~\ref{DDIp}; finally, the last equality recovers the description of the mapping $\mathrm{ip}^{(\llbracket 2\rrbracket,X)}$ introduced in Definition~\ref{DDVEch}.

Regarding the second item, let $s$ be a sort in $S$ and $P$ a term in $\mathrm{T}_{\Sigma}(X)_{s}$. The following chain of equalities holds
\allowdisplaybreaks
\begin{align*}
\mathrm{ip}^{(\llbracket 2\rrbracket,X)@}_{s}\left(
\eta^{(\llbracket 2\rrbracket,0)\sharp}_{s}\left(
P
\right)
\right)
&=
\mathrm{ip}^{(\llbracket 2\rrbracket,X)@}_{s}\left(
\left\llbracket
\eta^{( 2,0)\sharp}_{s}\left(
P
\right)
\right\rrbracket_{s}
\right)
\tag{1}
\\&=
\left\llbracket
\mathrm{ip}^{(2,X)@}_{s}\left(
\eta^{(2,0)\sharp}_{s}\left(
P
\right)
\right)
\right\rrbracket_{s}
\tag{2}
\\&=
\left\llbracket
\mathrm{ip}^{(2,0)\sharp}_{s}\left(
P
\right)
\right\rrbracket_{s}
\tag{3}
\\&=
\mathrm{ip}^{(\llbracket 2\rrbracket,0)\sharp}_{s}\left(
P
\right).
\tag{4}
\end{align*}

The first equality applies the mapping $\eta^{(\llbracket 2 \rrbracket,0)\sharp}$ introduced in Definition~\ref{DDPTQEtas}; the second equality applies the mapping $\mathrm{ip}^{(\llbracket 2\rrbracket,X)@}$ introduced in Definition~\ref{DDPTQIp}; the third equality follows from Proposition~\ref{PDIpDZ}; finally, the last equality recovers the description of the mapping $\mathrm{ip}^{(\llbracket 2\rrbracket,0)\sharp}$ introduced in Definition~\ref{DDVDZ}.

Regarding the third item, let $s$ be a sort in $S$ and $x$ a variable in $X_{s}$. The following chain of equalities holds
\allowdisplaybreaks
\begin{align*}
\mathrm{CH}^{\llbracket 2\rrbracket}_{s}\left(
\mathrm{ip}^{(\llbracket 2\rrbracket,X)}_{s}\left(
x
\right)
\right)
&=
\mathrm{CH}^{\llbracket 2\rrbracket}_{s}\left(
\bigl\llbracket
\mathrm{ip}^{(2,X)}_{s}(x)
\bigr\rrbracket_{s}
\right)
\tag{1}
\\&=
\biggl\llbracket
\mathrm{CH}^{(2)}_{s}\left(
\mathrm{ip}^{(2,X)}_{s}\left(
x
\right)
\right)
\biggr\rrbracket_{s}
\tag{2}
\\&=
\Bigl\llbracket
\eta^{(2,X)}_{s}\left(
x
\right)
\Bigr\rrbracket_{s}
\tag{3}
\\&=
\eta^{(\llbracket 2\rrbracket,X)}_{s}\left(
x
\right).
\tag{4}
\end{align*}

The first equality applies the mapping $\mathrm{ip}^{(\llbracket 2\rrbracket,X)}$ introduced in Definition~\ref{DDVEch}; the second equality applies the mapping $\mathrm{CH}^{\llbracket 2\rrbracket}$ introduced in Definition~\ref{DDPTQDCH}; the third equality follows from Proposition~\ref{PDCHDUId}; finally, the last equality recovers the mapping $\eta^{(\llbracket 2\rrbracket,X)}$ introduced in Definition~\ref{DDPTQEtas}.

Regarding the fourth item, let $s$ be a sort in $S$ and $P$ a term in $\mathrm{T}_{\Sigma}(X)_{s}$. The following chain of equalities holds
\allowdisplaybreaks
\begin{align*}
\mathrm{CH}^{\llbracket 2\rrbracket}_{s}\left(
\mathrm{ip}^{(\llbracket 2\rrbracket,0)\sharp}_{s}\left(
P
\right)
\right)
&=
\mathrm{CH}^{\llbracket 2\rrbracket}_{s}\left(
\left\llbracket
\mathrm{ip}^{(\llbracket 2\rrbracket,0)\sharp}_{s}\left(
P
\right)
\right\rrbracket_{s}
\right)
\tag{1}
\\&=
\left\llbracket
\mathrm{CH}^{(2)}_{s}\left(
\mathrm{ip}^{(2,0)\sharp}_{s}\left(
P
\right)
\right)
\right\rrbracket_{\Theta^{\llbracket 2\rrbracket}_{s}}
\tag{2}
\\&=
\left\llbracket
\eta^{(2,0)\sharp}_{s}\left(
P
\right)
\right\rrbracket_{\Theta^{\llbracket 2\rrbracket}_{s}}
\tag{3}
\\&=
\eta^{(\llbracket 2\rrbracket,0)\sharp}_{s}\left(
P
\right).
\tag{4}
\end{align*}

The first equality applies the mapping $\mathrm{ip}^{(\llbracket 2\rrbracket,0)\sharp}$ introduced in Definition~\ref{DDPTEtas}; the second equality applies the mapping $\mathrm{CH}^{\llbracket 2\rrbracket}$ introduced in Definition~\ref{DDPTQDCH}; the third equality follows from Proposition~\ref{PDCHDZId}; finally, the last equality recovers the mapping $\eta^{(\llbracket 2\rrbracket,0)\sharp}$ introduced in Definition~\ref{DDPTQEtas}.

This completes the proof.
\end{proof}

\begin{figure}
\begin{tikzpicture}
[ACliment/.style={-{To [angle'=45, length=5.75pt, width=4pt, round]}},scale=0.8]
\node[] (x) at (0,0) [] {$X$};
\node[] (t) at (6,0) [] {$\mathrm{T}_{\Sigma}(X)$};
\node[] (pth1) at (6,-2) [] {$[\mathrm{Pth}_{\boldsymbol{\mathcal{A}}}]$};
\node[] (pt1) at (6,-4) [] {$[\mathrm{PT}_{\boldsymbol{\mathcal{A}}}]$};

\node[] (pth2) at (6,-6) [] {$\llbracket \mathrm{Pth}_{\boldsymbol{\mathcal{A}}^{(2)}}\rrbracket$};
\node[] (pt2) at (6,-8) [] {$\llbracket\mathrm{PT}_{\boldsymbol{\mathcal{A}}^{(2)}}\rrbracket$};

\draw[ACliment]  (x) to node [above right]
{$\textstyle \eta^{(0,X)}$} (t);
\draw[ACliment, bend right=10]  (x) to node [pos=.7, fill=white]
{$\textstyle \mathrm{ip}^{([1],X)}$} (pth1.west);
\draw[ACliment, bend right=20]  (x) to node  [pos=.6, fill=white]
{$\textstyle\eta^{([1],X)}$} (pt1.west);

\draw[ACliment, bend right=30]  (x) to node [pos=.55, fill=white]
{$\textstyle \mathrm{ip}^{(\llbracket 2\rrbracket,X)}$} (pth2.west);
\draw[ACliment, bend right=40]  (x) to node  [below left]
{$\textstyle \eta^{(\llbracket 2\rrbracket,X)}$} (pt2.west);

\draw[ACliment] 
($(t)+(0,-.35)$) to node [right] {$\textstyle \mathrm{ip}^{([1],0)\sharp}$}  ($(pth1)+(0,.35)$);

\draw[ACliment] 
($(pth1)+(-.30,-.35)$) to node [left] {$\textstyle \mathrm{CH}^{[1]}$}  ($(pt1)+(-.30,.35)$);
\draw[ACliment] 
($(pt1)+(.30,+.35)$) to node [right] {$\textstyle \mathrm{ip}^{([1],X)@}$}  ($(pth1)+(.30,-.35)$);

\node[] () at ($(pt1)+(0,1)$) [] {$\cong$};

\draw[ACliment] 
($(pt1)+(0,-.35)$) to node [right] {$\textstyle \mathrm{ip}^{(\llbracket 2\rrbracket,[1])\sharp}$}  ($(pth2)+(0,.35)$);

\draw[ACliment] 
($(pth2)+(-.30,-.35)$) to node [left] {$\textstyle \mathrm{CH}^{\llbracket 2\rrbracket}$}  ($(pt2)+(-.30,.35)$);
\draw[ACliment] 
($(pt2)+(.30,+.35)$) to node [right] {$\textstyle \mathrm{ip}^{(\llbracket 2\rrbracket,X)@}$}  ($(pth2)+(.30,-.35)$);

\draw[ACliment, rounded corners]
(pt1.east) -- 
($(pt1)+(3,0)$) -- node [right] {$\textstyle \eta^{(\llbracket 2\rrbracket,[1])\sharp}$} 
($(pt2)+(3,+.1)$) -- ($(pt2.east)+(0,.1)$)
;

\draw[ACliment, rounded corners]
(t.east) -- 
($(t)+(5.5,0)$) -- node [right] {$\textstyle \eta^{(\llbracket 2\rrbracket,0)\sharp}$} 
($(pt2)+(5.5,-.1)$) -- 
($(pt2.east)-(0,.1)$)
;

\end{tikzpicture}
\caption{Quotient embeddings relative to $X$ at layers $0$, $1$ \& $2$.}
\label{FDPTQXEq}
\end{figure}

\begin{figure}
\begin{tikzpicture}
[ACliment/.style={-{To [angle'=45, length=5.75pt, width=4pt, round]}},scale=0.8]
\node[] (a) at (0,-2) [] {$\mathcal{A}$};
\node[] (pth1) at (6,-2) [] {$[\mathrm{Pth}_{\boldsymbol{\mathcal{A}}}]$};
\node[] (pt1) at (6,-4) [] {$[\mathrm{PT}_{\boldsymbol{\mathcal{A}}}]$};

\node[] (pth2) at (6,-6) [] {$\llbracket \mathrm{Pth}_{\boldsymbol{\mathcal{A}}^{(2)}}\rrbracket$};
\node[] (pt2) at (6,-8) [] {$\llbracket\mathrm{PT}_{\boldsymbol{\mathcal{A}}^{(2)}}\rrbracket$};

\draw[ACliment]  (a) to node [above right]
{$\textstyle \mathrm{ech}^{([1],\mathcal{A})}$} (pth1);
\draw[ACliment, bend right=10]  (a) to node [pos=.6, fill=white]
{$\textstyle \eta^{([1],\mathcal{A})}$} (pt1.west);
\draw[ACliment, bend right=20]  (a) to node  [midway, fill=white]
{$\textstyle \mathrm{ech}^{(\llbracket 2\rrbracket,\mathcal{A})}$} (pth2.west);

\draw[ACliment, bend right=30]  (a) to node [below left]
{$\textstyle \eta^{(\llbracket 2\rrbracket,\mathcal{A})}$} (pt2.west);

\draw[ACliment] 
($(pth1)+(-.30,-.35)$) to node [left] {$\textstyle \mathrm{CH}^{[1]}$}  ($(pt1)+(-.30,.35)$);
\draw[ACliment] 
($(pt1)+(.30,+.35)$) to node [right] {$\textstyle \mathrm{ip}^{([1],X)@}$}  ($(pth1)+(.30,-.35)$);

\node[] () at ($(pt1)+(0,1)$) [] {$\cong$};

\draw[ACliment] 
($(pt1)+(0,-.35)$) to node [right] {$\textstyle \mathrm{ip}^{(\llbracket 2\rrbracket,[1])\sharp}$}  ($(pth2)+(0,.35)$);

\draw[ACliment] 
($(pth2)+(-.30,-.35)$) to node [left] {$\textstyle \mathrm{CH}^{\llbracket 2\rrbracket}$}  ($(pt2)+(-.30,.35)$);
\draw[ACliment] 
($(pt2)+(.30,+.35)$) to node [right] {$\textstyle \mathrm{ip}^{(\llbracket 2\rrbracket,X)@}$}  ($(pth2)+(.30,-.35)$);

\draw[ACliment, rounded corners]
(pt1.east) -- 
($(pt1)+(3,0)$) -- node [right] {$\textstyle \eta^{(\llbracket 2\rrbracket,[1])\sharp}$} 
($(pt2)+(3,0)$) -- ($(pt2.east)+(0,0)$)
;

\end{tikzpicture}
\caption{Quotient embeddings relative to $\mathcal{A}$ at layers $1$ \& $2$.}
\label{FDPTQAEq}
\end{figure}

\begin{proposition}\label{PDPTQAEq} The following equalities holds
\begin{multicols}{2}
\begin{itemize}
\item[(i)] $\mathrm{ip}^{(\llbracket 2\rrbracket,X)@}\circ\eta^{(\llbracket 2\rrbracket,\mathcal{A})}=\mathrm{ech}^{(\llbracket 2\rrbracket,\mathcal{A})};$
\item[(ii)] $\mathrm{CH}^{\llbracket 2\rrbracket}\circ\mathrm{ech}^{(\llbracket 2\rrbracket,\mathcal{A})}=\eta^{(\llbracket 2\rrbracket,\mathcal{A})}$.
\end{itemize}
\end{multicols}
The reader is advised to consult the diagram presented in Figure~\ref{FDPTQAEq}.
\end{proposition}

\begin{proposition}\label{PDPTQDAEq} The following equalities holds
\begin{itemize}
\item[(i)] $\mathrm{ip}^{(\llbracket 2\rrbracket,X)@}\circ\eta^{(\llbracket 2\rrbracket,\mathcal{A}^{(2)})}=\mathrm{ech}^{(\llbracket 2\rrbracket,\mathcal{A}^{(2)})};$
\item[(ii)] $\mathrm{CH}^{\llbracket 2\rrbracket}\circ\mathrm{ech}^{(\llbracket 2\rrbracket,\mathcal{A}^{(2)})}=\eta^{(\llbracket 2\rrbracket,\mathcal{A}^{(2)})}$.
\end{itemize}
The reader is advised to consult the diagram presented in Figure~\ref{FDPTQDAEq}.
\end{proposition}

\begin{figure}
\begin{tikzpicture}
[ACliment/.style={-{To [angle'=45, length=5.75pt, width=4pt, round]}},scale=0.8]
\node[] (a) at (0,-2) [] {$\mathcal{A}^{(2)}$};
\node[] (pth1) at (6,-2) [] {$\llbracket \mathrm{Pth}_{\boldsymbol{\mathcal{A}}^{(2)}}\rrbracket$};
\node[] (pt1) at (6,-4) [] {$\llbracket\mathrm{PT}_{\boldsymbol{\mathcal{A}}^{(2)}}\rrbracket$};

\draw[ACliment]  (a) to node [above]
{$\textstyle \mathrm{ech}^{(\llbracket 2\rrbracket,\mathcal{A}^{(2)})}$} (pth1);
\draw[ACliment, bend right=10]  (a) to node  [below left]
{$\textstyle \eta^{(\llbracket 2\rrbracket,\mathcal{A}^{(2)})}$} (pt1.west);

\draw[ACliment] 
($(pth1)+(-.3,-.35)$) to node [left] {$\textstyle \mathrm{CH}^{\llbracket 2\rrbracket}$}  ($(pt1)+(-.3,.35)$);
\draw[ACliment] 
($(pt1)+(.3,+.35)$) to node [right] {$\textstyle \mathrm{ip}^{(\llbracket 2\rrbracket,X)@}$}  ($(pth1)+(.3,-.35)$);
\end{tikzpicture}
\caption{Quotient embeddings relative to $\mathcal{A}^{(2)}$ at layer $2$.}
\label{FDPTQDAEq}
\end{figure}

\section{
\texorpdfstring
{A structure of partial $\Sigma^{\boldsymbol{\mathcal{A}}^{(2)}}$-algebra on $\llbracket\mathrm{PT}_{\boldsymbol{\mathcal{A}}^{(2)}}\rrbracket$}
{A  partial algebra on the quotient of second-order path terms}
}

In this section we study  the structure of partial $\Sigma^{\boldsymbol{\mathcal{A}}^{(2)}}$-algebra structure that naturally arises in the many-sorted set of second-order path term classes, which is a subalgebra of $\mathbf{T}_{\Sigma^{\boldsymbol{\mathcal{A}}^{(2)}}}(X)/{\Theta^{\llbracket 2 \rrbracket}}$. 

\begin{restatable}{proposition}{PDPTQDCatAlg}
\label{PDPTQDCatAlg} 
\index{path terms!second-order!$\llbracket\mathbf{PT}_{\boldsymbol{\mathcal{A}}^{(2)}}\rrbracket$}
The $S$-sorted set $\llbracket\mathrm{PT}_{\boldsymbol{\mathcal{A}}^{(2)}}\rrbracket$ is equipped, in a natural way, with a structure of many-sorted partial $\Sigma^{\boldsymbol{\mathcal{A}}^{(2)}}$-algebra, which is a $\Sigma^{\boldsymbol{\mathcal{A}}^{(2)}}$-subalgebra of $\mathbf{T}_{\Sigma^{\boldsymbol{\mathcal{A}}^{(2)}}}(X)/{\Theta^{\llbracket 2 \rrbracket}}$.
\end{restatable}
\begin{proof}
Let us denote by $\llbracket\mathbf{PT}_{\boldsymbol{\mathcal{A}}^{(2)}}\rrbracket$ the partial $\Sigma^{\boldsymbol{\mathcal{A}}^{(2)}}$-algebra defined as follows:

\textsf{(i)} The underlying $S$-sorted set of $\llbracket\mathbf{PT}_{\boldsymbol{\mathcal{A}}^{(2)}}\rrbracket$ is $\llbracket\mathrm{PT}_{\boldsymbol{\mathcal{A}}^{(2)}}\rrbracket$.

\textsf{(ii)} For every $(\mathbf{s},s)\in S^{\star}\times S$ and every operation symbol $\sigma\in\Sigma_{\mathbf{s},s}$, the operation $\sigma^{\llbracket\mathbf{PT}_{\boldsymbol{\mathcal{A}}^{(2)}}\rrbracket}$ is equal to $\sigma^{\mathbf{T}_{\Sigma^{\boldsymbol{\mathcal{A}}^{(2)}}}(X)/{\Theta^{\llbracket 2 \rrbracket}}}$, i.e., to the interpretation of $\sigma$ in $\mathbf{T}_{\Sigma^{\boldsymbol{\mathcal{A}}^{(2)}}}(X)/{\Theta^{\llbracket 2 \rrbracket}}$. 

Let us note that, if $(\llbracket P_{j}\rrbracket_{s_{j}})_{j\in\bb{\mathbf{s}}}$ is a family of second-order path term classes in $\llbracket \mathrm{PT}_{\boldsymbol{\mathcal{A}}^{(2)}}\rrbracket_{\mathbf{s}}$, then the term class
$$
\sigma^{\mathbf{T}_{\Sigma^{\boldsymbol{\mathcal{A}}^{(2)}}}(X)/{\Theta^{\llbracket 2 \rrbracket}}}
\left(\left(
\llbracket P_{j}
\rrbracket_{s_{j}}
\right)_{j\in\bb{\mathbf{s}}}\right)
=
\biggl\llbracket
\sigma^{\mathbf{T}_{\Sigma^{\boldsymbol{\mathcal{A}}^{(2)}}}(X)}
\left(\left(P_{j}
\right)_{j\in\bb{\mathbf{s}}}\right)
\biggr\rrbracket_{s}
$$
is a second-order path term class by Proposition~\ref{PDPTDCatAlg}, i.e.,  $\sigma^{\llbracket\mathbf{PT}_{\boldsymbol{\mathcal{A}}^{(2)}}\rrbracket}$ is well-defined.

\textsf{(iii)} For every $s\in S$ and every rewrite rule $\mathfrak{p}\in \mathcal{A}_{s}$, the constant $\mathfrak{p}^{\llbracket\mathbf{PT}_{\boldsymbol{\mathcal{A}}^{(2)}}\rrbracket}$ is equal to $\mathfrak{p}^{\mathbf{T}_{\Sigma^{\boldsymbol{\mathcal{A}}^{(2)}}}(X)/{\Theta^{\llbracket 2 \rrbracket}}}$, i.e., to the interpretation of $\mathfrak{p}$ in $\mathbf{T}_{\Sigma^{\boldsymbol{\mathcal{A}}^{(2)}}}(X)/{\Theta^{\llbracket 2 \rrbracket}}$. 

Let us note that the term class
$$
\mathfrak{p}^{\llbracket\mathbf{PT}_{\boldsymbol{\mathcal{A}}^{(2)}}\rrbracket}
=
\Bigl\llbracket
\mathfrak{p}^{\mathbf{T}_{\Sigma^{\boldsymbol{\mathcal{A}}^{(2)}}}(X)}
\Bigr\rrbracket_{s}
$$
is a second-order path term class by Proposition~\ref{PDPTDCatAlg}, i.e., $
\mathfrak{p}^{\llbracket\mathbf{PT}_{\boldsymbol{\mathcal{A}}^{(2)}}\rrbracket}$ is well-defined.

\textsf{(iv)} For every $s\in S$, the operation $\mathrm{sc}_{s}^{0\llbracket\mathbf{PT}_{\boldsymbol{\mathcal{A}}^{(2)}}\rrbracket}$ is equal to $\mathrm{sc}_{s}^{0\mathbf{T}_{\Sigma^{\boldsymbol{\mathcal{A}}^{(2)}}}(X)/{\Theta^{\llbracket 2 \rrbracket}}}$, i.e., to the interpretation of $\mathrm{sc}_{s}^{0}$ in $\mathbf{T}_{\Sigma^{\boldsymbol{\mathcal{A}}^{(2)}}}(X)/{\Theta^{\llbracket 2 \rrbracket}}$. 

Let us note that, if $\llbracket P\rrbracket_{s}$ is a second-order path term class in $\llbracket \mathrm{PT}_{\boldsymbol{\mathcal{A}}^{(2)}}\rrbracket_{s}$, then the term class
$$
\mathrm{sc}_{s}^{0\mathbf{T}_{\Sigma^{\boldsymbol{\mathcal{A}}^{(2)}}}(X)/{\Theta^{\llbracket 2 \rrbracket}}}
\left(
\llbracket P
\rrbracket_{s}
\right)
=
\biggl\llbracket
\mathrm{sc}_{s}^{0\mathbf{T}_{\Sigma^{\boldsymbol{\mathcal{A}}^{(2)}}}(X)}
\left(P
\right)
\biggr\rrbracket_{s}
$$
is a second-order path term class by Proposition~\ref{PDPTDCatAlg}. That is, the operation $\mathrm{sc}_{s}^{0\llbracket\mathbf{PT}_{\boldsymbol{\mathcal{A}}^{(2)}}\rrbracket}$ is well-defined.

\textsf{(v)} For every $s\in S$, the operation $\mathrm{tg}_{s}^{0\llbracket\mathbf{PT}_{\boldsymbol{\mathcal{A}}^{(2)}}\rrbracket}$ is equal to $\mathrm{tg}_{s}^{0\mathbf{T}_{\Sigma^{\boldsymbol{\mathcal{A}}^{(2)}}}(X)/{\Theta^{\llbracket 2 \rrbracket}}}$, i.e., to the interpretation of $\mathrm{tg}_{s}^{0}$ in $\mathbf{T}_{\Sigma^{\boldsymbol{\mathcal{A}}^{(2)}}}(X)/{\Theta^{\llbracket 2 \rrbracket}}$.  Thus $\mathrm{tg}_{s}^{0\llbracket\mathbf{PT}_{\boldsymbol{\mathcal{A}}^{(2)}}\rrbracket}$ is well-defined by a similar argument to that of $\mathrm{sc}_{s}^{0\llbracket\mathbf{PT}_{\boldsymbol{\mathcal{A}}^{(2)}}\rrbracket}$.

\textsf{(vi)} For every $s\in S$, the partial binary operation $\circ_{s}^{0\llbracket\mathbf{PT}_{\boldsymbol{\mathcal{A}}^{(2)}}\rrbracket}$ is defined on second-order path term classes $\llbracket Q\rrbracket_{s}$, $\llbracket P\rrbracket_{s}$ in $\llbracket\mathrm{PT}_{\boldsymbol{\mathcal{A}}^{(2)}}\rrbracket_{s}$ if, and only if, the following equation holds
$$
\mathrm{sc}_{s}^{(0,2)}\left(\mathrm{ip}^{(2,X)@}_{s}(Q)\right)=
\mathrm{tg}_{s}^{(0,2)}\left(\mathrm{ip}^{(2,X)@}_{s}(P)\right).
$$

Note that the last equation does not depend on the representative in the $\Theta^{\llbracket 2\rrbracket}_{s}$-class by Corollary~\ref{CDPTScTg}.

For this case, $\circ_{s}^{0 \llbracket\mathbf{PT}_{\boldsymbol{\mathcal{A}}^{(2)}}\rrbracket}$ is equal to $\circ_{s}^{0\mathbf{T}_{\Sigma^{\boldsymbol{\mathcal{A}}^{(2)}}}(X)/{\Theta^{\llbracket 2 \rrbracket}}}$, i.e., to the interpretation of $\circ^{0}_{s}$ in $\mathbf{T}_{\Sigma^{\boldsymbol{\mathcal{A}}^{(2)}}}(X)/{\Theta^{\llbracket 2 \rrbracket}}$.  Hence, if $\llbracket Q\rrbracket_{s}$, $\llbracket P\rrbracket_{s}$ are second-order path term classes  in $\llbracket\mathrm{PT}_{\boldsymbol{\mathcal{A}}^{(2)}}\rrbracket_{s}$ with
$
\mathrm{sc}_{s}^{(0,2)}(\mathrm{ip}^{(2,X)@}_{s}(Q))=
\mathrm{tg}_{s}^{(0,2)}(\mathrm{ip}^{(2,X)@}_{s}(P)),
$ 
then the term class
$$
\llbracket Q\rrbracket_{s}
\circ_{s}^{0\mathbf{T}_{\Sigma^{\boldsymbol{\mathcal{A}}^{(2)}}}(X)/{\Theta^{\llbracket 2 \rrbracket}}}
\llbracket P\rrbracket_{s}
=
\biggl\llbracket
Q
\circ_{s}^{0\mathbf{T}_{\Sigma^{\boldsymbol{\mathcal{A}}^{(2)}}}(X)}
P
\biggr\rrbracket_{s}
$$
is a second-order path term class by Proposition~\ref{PDPTDCatAlg}, i.e., $\circ_{s}^{0\llbracket\mathbf{PT}_{\boldsymbol{\mathcal{A}}^{(2)}}\rrbracket}$ is well-defined.

\textsf{(vii)} For every $s\in S$ and every second-order rewrite rule $\mathfrak{p}^{(2)}\in \mathcal{A}^{(2)}_{s}$, the constant $\mathfrak{p}^{(2)\llbracket\mathbf{PT}_{\boldsymbol{\mathcal{A}}^{(2)}}\rrbracket}$ is equal to $\mathfrak{p}^{(2)\mathbf{T}_{\Sigma^{\boldsymbol{\mathcal{A}}^{(2)}}}(X)/{\Theta^{\llbracket 2 \rrbracket}}}$, i.e., to the interpretation of $\mathfrak{p}^{(2)}$ in $\mathbf{T}_{\Sigma^{\boldsymbol{\mathcal{A}}^{(2)}}}(X)/{\Theta^{\llbracket 2 \rrbracket}}$. 

Let us note that, by Proposition~\ref{PDPTDCatAlg}, the term class
$$
\mathfrak{p}^{(2)\llbracket\mathbf{PT}_{\boldsymbol{\mathcal{A}}^{(2)}}\rrbracket}
=
\Bigl\llbracket
\mathfrak{p}^{(2)\mathbf{T}_{\Sigma^{\boldsymbol{\mathcal{A}}^{(2)}}}(X)}
\Bigr\rrbracket_{s}
$$
is a second-order path term class. That is, the constant $
\mathfrak{p}^{(2)\llbracket\mathbf{PT}_{\boldsymbol{\mathcal{A}}^{(2)}}\rrbracket}$ is well-defined.

\textsf{(viii)} For every $s\in S$, the operation $\mathrm{sc}_{s}^{1\llbracket\mathbf{PT}_{\boldsymbol{\mathcal{A}}^{(2)}}\rrbracket}$ is equal to $\mathrm{sc}_{s}^{1\mathbf{T}_{\Sigma^{\boldsymbol{\mathcal{A}}^{(2)}}}(X)/{\Theta^{\llbracket 2 \rrbracket}}}$, i.e., to the interpretation of $\mathrm{sc}_{s}^{1}$ in $\mathbf{T}_{\Sigma^{\boldsymbol{\mathcal{A}}^{(2)}}}(X)/{\Theta^{\llbracket 2 \rrbracket}}$. 

Let us note that, if $\llbracket P\rrbracket_{s}$ is a  second-order path term class in $\llbracket\mathrm{PT}_{\boldsymbol{\mathcal{A}}^{(2)}}\rrbracket_{s}$, then the term class
$$
\mathrm{sc}_{s}^{1\mathbf{T}_{\Sigma^{\boldsymbol{\mathcal{A}}^{(2)}}}(X)/{\Theta^{\llbracket 2 \rrbracket}}}
\left(
\llbracket P\rrbracket_{s}
\right)
=
\biggl\llbracket 
\mathrm{sc}_{s}^{1\mathbf{T}_{\Sigma^{\boldsymbol{\mathcal{A}}^{(2)}}}(X)}
(P)
\biggr\rrbracket_{s}
$$
is a second-order path term class by Proposition~\ref{PDPTDCatAlg}. That is, the operation $\mathrm{sc}_{s}^{1\llbracket\mathbf{PT}_{\boldsymbol{\mathcal{A}}^{(2)}}\rrbracket}$ is well-defined.

\textsf{(ix)} For every $s\in S$, the operation $\mathrm{tg}_{s}^{1\llbracket\mathbf{PT}_{\boldsymbol{\mathcal{A}}^{(2)}}\rrbracket}$ is equal to $\mathrm{tg}_{s}^{1\mathbf{T}_{\Sigma^{\boldsymbol{\mathcal{A}}^{(2)}}}(X)/{\Theta^{\llbracket 2 \rrbracket}}}$, i.e., to the interpretation of $\mathrm{tg}_{s}^{1}$ in $\mathbf{T}_{\Sigma^{\boldsymbol{\mathcal{A}}^{(2)}}}(X)/{\Theta^{\llbracket 2 \rrbracket}}$.  Thus $\mathrm{tg}_{s}^{1\llbracket\mathbf{PT}_{\boldsymbol{\mathcal{A}}^{(2)}}\rrbracket}$ is well-defined by a similar argument to that of $\mathrm{sc}_{s}^{1\llbracket\mathbf{PT}_{\boldsymbol{\mathcal{A}}^{(2)}}\rrbracket}$.

\textsf{(x)} For every $s\in S$, the partial binary operation $\circ_{s}^{1\llbracket\mathbf{PT}_{\boldsymbol{\mathcal{A}}^{(2)}}\rrbracket}$ is defined on second-order  path term classes $\llbracket Q\rrbracket_{s}$, $\llbracket P\rrbracket_{s}$ in $\llbracket\mathrm{PT}_{\boldsymbol{\mathcal{A}}^{(2)}}\rrbracket_{s}$ if, and only if, the following equation holds
$$
\mathrm{sc}_{s}^{(1,2)}\left(\mathrm{ip}^{(2,X)@}_{s}(Q)\right)=
\mathrm{tg}_{s}^{(1,2)}\left(\mathrm{ip}^{(2,X)@}_{s}(P)\right).
$$

Note that the last equation does not depend on the representative in the $\Theta^{\llbracket 2 \rrbracket}_{s}$-class by Corollary~\ref{CDPTScTg}.

For this case, $\circ_{s}^{1\llbracket \mathbf{PT}_{\boldsymbol{\mathcal{A}}^{(2)}}\rrbracket}$ is equal to $\circ_{s}^{1\mathbf{T}_{\Sigma^{\boldsymbol{\mathcal{A}}^{(2)}}}(X)/{\Theta^{\llbracket 2 \rrbracket}}}$, i.e., to the interpretation of $\circ^{1}_{s}$ in $\mathbf{T}_{\Sigma^{\boldsymbol{\mathcal{A}}^{(2)}}}(X)/{\Theta^{\llbracket 2 \rrbracket}}$.  Hence, if $\llbracket Q\rrbracket_{s}$, $\llbracket P\rrbracket_{s}$ are  second-order path term classes in $\llbracket\mathrm{PT}_{\boldsymbol{\mathcal{A}}^{(2)}}\rrbracket_{s}$ with
$
\mathrm{sc}_{s}^{(1,2)}(\mathrm{ip}^{(2,X)@}_{s}(Q))=
\mathrm{tg}_{s}^{(1,2)}(\mathrm{ip}^{(2,X)@}_{s}(P)),
$ 
then the term class
$$
\llbracket Q\rrbracket_{s}
\circ_{s}^{1\mathbf{T}_{\Sigma^{\boldsymbol{\mathcal{A}}^{(2)}}}(X)/{\Theta^{\llbracket 2 \rrbracket}}}
\llbracket P\rrbracket_{s}
=
\biggl\llbracket
Q
\circ_{s}^{1\mathbf{T}_{\Sigma^{\boldsymbol{\mathcal{A}}^{(2)}}}(X)}
P
\biggr\rrbracket_{s}
$$
is a second-order path term class by Proposition~\ref{PDPTDCatAlg}. That is, the operation $\circ_{s}^{1\llbracket\mathbf{PT}_{\boldsymbol{\mathcal{A}}^{(2)}}\rrbracket}$ is well-defined.

This completes the definition of $\llbracket\mathbf{PT}_{\boldsymbol{\mathcal{A}}^{(2)}}\rrbracket$, the partial many-sorted $\Sigma^{\boldsymbol{\mathcal{A}}^{(2)}}$-algebra of second-order path term classes.
\end{proof}

We next state that for the many-sorted partial $\Sigma^{\boldsymbol{\mathcal{A}}^{(2)}}$-algebras $\mathbf{PT}_{\boldsymbol{\mathcal{A}}^{(2)}}$, of second-order path terms, and $\llbracket\mathbf{PT}_{\boldsymbol{\mathcal{A}}^{(2)}}\rrbracket$, of second-order path term classes, introduced in Proposition~\ref{PDPTDCatAlg} and Proposition~\ref{PDPTQDCtyAlg}, respectively, the canonical projection for the $\Sigma^{\boldsymbol{\mathcal{A}}^{(2)}}$-congruence $\Theta^{\llbracket 2 \rrbracket}$ determines a surjective $\Sigma^{\boldsymbol{\mathcal{A}}^{(2)}}$-homomorphism. 

\begin{restatable}{corollary}{CDPTQPr}
\label{CDPTQPr} The mapping $\mathrm{pr}^{\Theta^{\llbracket 2 \rrbracket}}$ is a surjective $\Sigma^{\boldsymbol{\mathcal{A}}^{(2)}}$-homomorphism
$$
\mathrm{pr}^{\Theta^{\llbracket 2 \rrbracket}}\colon
\mathbf{PT}_{\boldsymbol{\mathcal{A}}^{(2)}}
\mor
\llbracket\mathbf{PT}_{\boldsymbol{\mathcal{A}}^{(2)}}\rrbracket.
$$
\end{restatable}

\section{
\texorpdfstring
{A structure of $2$-categorial $\Sigma$-algebra on $\llbracket\mathrm{PT}_{\boldsymbol{\mathcal{A}}^{(2)}}\rrbracket$}
{A categorial algebra on the quotient of second-order path terms}
}

In this section we study the algebraic and categorial structure of the second-order path term quotient by $\Theta^{\llbracket 2\rrbracket}$. We begin by presenting the equations that satisfy the partial many-sorted $\Sigma^{\boldsymbol{\mathcal{A}}^{(2)}}$-algebra $\llbracket\mathbf{PT}_{\boldsymbol{\mathcal{A}}^{(2)}}\rrbracket$.

\begin{proposition}\label{PDPTQVarA2} Let $s$ be a sort in $S$ and $\llbracket P\rrbracket_{s}$ a second-order path term class in  $\llbracket\mathrm{PT}_{\boldsymbol{\mathcal{A}}^{(2)}}\rrbracket_{s}$, then the following equalities hold
\allowdisplaybreaks
\begin{align*}
\mathrm{sc}^{0\llbracket\mathbf{PT}_{\boldsymbol{\mathcal{A}}^{(2)}}\rrbracket}_{s}\left(
\mathrm{sc}^{0\llbracket\mathbf{PT}_{\boldsymbol{\mathcal{A}}^{(2)}}\rrbracket}_{s}\left(
\llbracket P\rrbracket_{s}
\right)
\right)
& 
=
\mathrm{sc}^{0\llbracket\mathbf{PT}_{\boldsymbol{\mathcal{A}}^{(2)}}\rrbracket}_{s}\left(
\llbracket P\rrbracket_{s}
\right);
\\
\mathrm{sc}^{0\llbracket\mathbf{PT}_{\boldsymbol{\mathcal{A}}^{(2)}}\rrbracket}_{s}\left(
\mathrm{tg}^{0\llbracket\mathbf{PT}_{\boldsymbol{\mathcal{A}}^{(2)}}\rrbracket}_{s}\left(
\llbracket P\rrbracket_{s}
\right)
\right)
& 
=
\mathrm{tg}^{0\llbracket\mathbf{PT}_{\boldsymbol{\mathcal{A}}^{(2)}}\rrbracket}_{s}\left(
\llbracket P\rrbracket_{s}
\right);
\\
\mathrm{tg}^{0\llbracket\mathbf{PT}_{\boldsymbol{\mathcal{A}}^{(2)}}\rrbracket}_{s}\left(
\mathrm{sc}^{0\llbracket\mathbf{PT}_{\boldsymbol{\mathcal{A}}^{(2)}}\rrbracket}_{s}\left(
\llbracket P\rrbracket_{s}
\right)
\right)
& 
=
\mathrm{sc}^{0\llbracket\mathbf{PT}_{\boldsymbol{\mathcal{A}}^{(2)}}\rrbracket}_{s}\left(
\llbracket P\rrbracket_{s}
\right);
\\
\mathrm{tg}^{0\llbracket\mathbf{PT}_{\boldsymbol{\mathcal{A}}^{(2)}}\rrbracket}_{s}\left(
\mathrm{tg}^{0\llbracket\mathbf{PT}_{\boldsymbol{\mathcal{A}}^{(2)}}\rrbracket}_{s}\left(
\llbracket P\rrbracket_{s}
\right)
\right)
& 
=
\mathrm{tg}^{0\llbracket\mathbf{PT}_{\boldsymbol{\mathcal{A}}^{(2)}}\rrbracket}_{s}\left(
\llbracket P\rrbracket_{s}
\right).
\end{align*}	
\end{proposition}
\begin{proof}
We will only present the proof for the first equality. The other three are handled in a similar manner.

Note that the following chain of equalities holds
\begin{flushleft}
$\mathrm{sc}^{0\llbracket\mathbf{PT}_{\boldsymbol{\mathcal{A}}^{(2)}}\rrbracket}_{s}\left(
\mathrm{sc}^{0\llbracket\mathbf{PT}_{\boldsymbol{\mathcal{A}}^{(2)}}\rrbracket}_{s}\left(
\llbracket P\rrbracket_{s}
\right)
\right)$
\allowdisplaybreaks
\begin{align*}
\quad&=
\biggl\llbracket
\mathrm{sc}^{0\mathbf{PT}_{\boldsymbol{\mathcal{A}}^{(2)}}}_{s}\left(
\mathrm{sc}^{0\mathbf{PT}_{\boldsymbol{\mathcal{A}}^{(2)}}}_{s}\left(
P
\right)\right)
\biggr\rrbracket_{s}
\tag{1}
\\&=
\biggl\llbracket
\mathrm{CH}^{(2)}_{s}\left(
\mathrm{ip}^{(2,X)@}_{s}\left(
\mathrm{sc}^{0\mathbf{PT}_{\boldsymbol{\mathcal{A}}^{(2)}}}_{s}\left(
\mathrm{sc}^{0\mathbf{PT}_{\boldsymbol{\mathcal{A}}^{(2)}}}_{s}\left(
P
\right)\right)
\right)\right)
\biggr\rrbracket_{s}
\tag{2}
\\&=
\biggl\llbracket
\mathrm{CH}^{(2)}_{s}\left(
\mathrm{sc}^{0\mathbf{Pth}_{\boldsymbol{\mathcal{A}}^{(2)}}}_{s}\left(
\mathrm{sc}^{0\mathbf{Pth}_{\boldsymbol{\mathcal{A}}^{(2)}}}_{s}\left(
\mathrm{ip}^{(2,X)@}_{s}\left(
P
\right)\right)
\right)\right)
\biggr\rrbracket_{s}
\tag{3}
\\&=
\biggl\llbracket
\mathrm{CH}^{(2)}_{s}\left(
\mathrm{sc}^{0\mathbf{Pth}_{\boldsymbol{\mathcal{A}}^{(2)}}}_{s}\left(
\mathrm{ip}^{(2,X)@}_{s}\left(
P
\right)
\right)\right)
\biggr\rrbracket_{s}
\tag{4}
\\&=
\biggl\llbracket
\mathrm{CH}^{(2)}_{s}\left(
\mathrm{ip}^{(2,X)@}_{s}\left(
\mathrm{sc}^{0\mathbf{PT}_{\boldsymbol{\mathcal{A}}^{(2)}}}_{s}\left(
P
\right)
\right)\right)
\biggr\rrbracket_{s}
\tag{5}
\\&=
\biggl\llbracket
\mathrm{sc}^{0\mathbf{PT}_{\boldsymbol{\mathcal{A}}^{(2)}}}_{s}\left(
P
\right)
\biggr\rrbracket_{s}
\tag{6}
\\&=
\mathrm{sc}^{0\llbracket\mathbf{PT}_{\boldsymbol{\mathcal{A}}^{(2)}}\rrbracket}_{s}\left(
\llbracket P\rrbracket_{s}
\right).
\tag{7}
\end{align*}
\end{flushleft}

In the just stated chain of equalities, the first equality unravels the description of the operation of $0$-source in the many sorted partial $\Sigma^{\boldsymbol{\mathcal{A}}^{(2)}}$-algebra $\llbracket\mathbf{PT}_{\boldsymbol{\mathcal{A}}^{(2)}}\rrbracket$ according to Proposition~\ref{PDPTQDCatAlg}; the second equality follows from Proposition~\ref{LDVWCong}; the third equality follows from the fact that $\mathrm{ip}^{(2,X)@}$ is a $\Sigma^{\boldsymbol{\mathcal{A}}^{(2)}}$-homomorphism according to Definition~\ref{DDIp}; the fourth equality follows from Proposition~\ref{PDVDCH} and Proposition~\ref{PDVVarA2}; the fifth equality follows from the fact that $\mathrm{ip}^{(2,X)@}$ is a $\Sigma^{\boldsymbol{\mathcal{A}}^{(2)}}$-homomorphism according to Definition~\ref{DDIp}; the sixth equality follows from Proposition~\ref{LDVWCong}; finally, the last equality recovers the description of the operation of $0$-source in the many sorted partial $\Sigma^{\boldsymbol{\mathcal{A}}^{(2)}}$-algebra $\llbracket\mathbf{PT}_{\boldsymbol{\mathcal{A}}^{(2)}}\rrbracket$ according to Proposition~\ref{PDPTQDCatAlg}.

This completes the proof.
\end{proof}

\begin{proposition}\label{PDPTQVarA3} Let $s$ be a sort in $S$ and $\llbracket P\rrbracket_{s}$, $\llbracket Q\rrbracket_{s}$ second-order path term classes in  $\llbracket\mathrm{PT}_{\boldsymbol{\mathcal{A}}^{(2)}}\rrbracket_{s}$, then the following statement are equivalent
\begin{itemize}
\item[(i)] $\llbracket Q\rrbracket_{s}
\circ^{0\llbracket\mathbf{PT}_{\boldsymbol{\mathcal{A}}^{(2)}}\rrbracket}_{s}
\llbracket P\rrbracket_{s}
$ is defined;
\item[(ii)] $\mathrm{sc}^{0\llbracket\mathbf{PT}_{\boldsymbol{\mathcal{A}}^{(2)}}\rrbracket}_{s}(
\llbracket Q\rrbracket_{s}
)=
\mathrm{tg}^{0\llbracket\mathbf{PT}_{\boldsymbol{\mathcal{A}}^{(2)}}\rrbracket}_{s}(
\llbracket P\rrbracket_{s}
)
$.
\end{itemize}
\end{proposition}
\begin{proof}
The following chain of equivalences holds
\begin{flushleft}
$\llbracket Q\rrbracket_{s}
\circ^{0\llbracket\mathbf{PT}_{\boldsymbol{\mathcal{A}}^{(2)}}\rrbracket}_{s}
\llbracket P\rrbracket_{s}$
\mbox{ is defined}
\allowdisplaybreaks
\begin{align*}
&\Leftrightarrow
\mathrm{sc}^{(0,2)}_{s}\left(\mathrm{ip}^{(2,X)@}_{s}\left(
Q\right)\right)
=
\mathrm{tg}^{(0,2)}_{s}\left(\mathrm{ip}^{(2,X)@}_{s}\left(
P\right)\right)
\tag{1}
\\&\Leftrightarrow
\mathrm{sc}^{0\llbracket \mathbf{Pth}_{\boldsymbol{\mathcal{A}}^{(2)}} \rrbracket}_{s}\left(
\Bigl\llbracket \mathrm{ip}^{(2,X)@}_{s}\left(Q\right) \Bigr\rrbracket_{s}\right)
=
\mathrm{tg}^{0\llbracket \mathbf{Pth}_{\boldsymbol{\mathcal{A}}^{(2)}} \rrbracket}_{s}\left(
\Bigl\llbracket \mathrm{ip}^{(2,X)@}_{s}\left(P\right) \Bigr\rrbracket_{s}\right)
\tag{2}
\\&\Leftrightarrow
\Bigl\llbracket \mathrm{sc}^{0\mathbf{Pth}_{\boldsymbol{\mathcal{A}}^{(2)}}}_{s}\left(\mathrm{ip}^{(2,X)@}_{s}\left(Q\right)\right) \Bigr\rrbracket_{s}
=
\Bigl\llbracket \mathrm{tg}^{0\mathbf{Pth}_{\boldsymbol{\mathcal{A}}^{(2)}}}_{s}\left(\mathrm{ip}^{(2,X)@}_{s}\left(P\right)\right) \Bigr\rrbracket_{s}
\tag{3}
\\&\Leftrightarrow
\biggl\llbracket \mathrm{ip}^{(2,X)@}_{s}\left( \mathrm{sc}^{0\mathbf{PT}_{\boldsymbol{\mathcal{A}}^{(2)}}}_{s}\left(Q\right)\right) \biggr\rrbracket_{s}=
\biggl\llbracket \mathrm{ip}^{(2,X)@}_{s}\left( \mathrm{tg}^{0\mathbf{PT}_{\boldsymbol{\mathcal{A}}^{(2)}}}_{s}\left(P\right)\right) \biggr\rrbracket_{s}
\tag{4}
\\&\Leftrightarrow
\biggl\llbracket \mathrm{CH}^{(2)}_{s}\left(\mathrm{ip}^{(2,X)@}_{s}\left( \mathrm{sc}^{0\mathbf{PT}_{\boldsymbol{\mathcal{A}}^{(2)}}}_{s}\left(Q\right)\right)\right) \biggr\rrbracket_{s}
\\&\qquad\qquad\qquad\qquad=
\biggl\llbracket \mathrm{CH}^{(2)}_{s}\left(\mathrm{ip}^{(2,X)@}_{s}\left( \mathrm{tg}^{0\mathbf{PT}_{\boldsymbol{\mathcal{A}}^{(2)}}}_{s}\left(P\right)\right)\right) \biggr\rrbracket_{s}
\tag{5}
\\&\Leftrightarrow
\Bigl\llbracket  \mathrm{sc}^{0\mathbf{PT}_{\boldsymbol{\mathcal{A}}^{(2)}}}_{s}\left(Q\right) \Bigr\rrbracket_{s}
=
\Bigl\llbracket  \mathrm{tg}^{0\mathbf{PT}_{\boldsymbol{\mathcal{A}}^{(2)}}}_{s}\left(P\right) \Bigr\rrbracket_{s}
\tag{6}
\\&\Leftrightarrow
\mathrm{sc}^{0\llbracket\mathbf{PT}_{\boldsymbol{\mathcal{A}}^{(2)}}\rrbracket}_{s}\left(
\llbracket Q\rrbracket_{s}
\right)=
\mathrm{tg}^{0\llbracket\mathbf{PT}_{\boldsymbol{\mathcal{A}}^{(2)}}\rrbracket}_{s}\left(
\llbracket P\rrbracket_{s}
\right).
\tag{7}
\end{align*}
\end{flushleft}

In the just stated chain of equivalences, the first equivalence follows from the description of the $0$-composition operation in the many-sorted partial $\Sigma^{\boldsymbol{\mathcal{A}}^{(2)}}$-algebra $\llbracket\mathbf{PT}_{\boldsymbol{\mathcal{A}}^{(2)}}\rrbracket$ according to Proposition~\ref{PDPTQDCatAlg}; the second equivalence follows from Proposition~\ref{PDVVarA3}; the third equivalence simply unravels the description of the $0$-source and $0$-target operations in the many-sorted partial $\Sigma^{\boldsymbol{\mathcal{A}}^{(2)}}$-algebra $\llbracket \mathbf{Pth}_{\boldsymbol{\mathcal{A}}^{(2)}} \rrbracket$, according to Proposition~\ref{PDVDCatAlg}; the fourth equivalence follows from the fact that, according to Definition~\ref{DDIp}, $\mathrm{ip}^{(2,X)@}$ is a $\Sigma^{\boldsymbol{\mathcal{A}}^{(2)}}$-homomorphism; the fifth equivalence follows from left to right from Proposition~\ref{PDVDCH} and from right to left by Corollary~\ref{CDVCongDCH}; the sixth equivalence follows from Lemma~\ref{LDVWCong}; finally, the last equivalence recovers the description of the $0$-source and $0$-target operations in the many-sorted partial $\Sigma^{\boldsymbol{\mathcal{A}}^{(2)}}$-algebra $\llbracket\mathbf{PT}_{\boldsymbol{\mathcal{A}}^{(2)}}\rrbracket$ according to Proposition~\ref{PDPTQDCatAlg}.

This completes the proof.
\end{proof}

\begin{proposition}\label{PDPTQVarA4} Let $s$ be a sort in $S$ and $\llbracket P\rrbracket_{s}$, $\llbracket Q\rrbracket_{s}$ second-order path term classes in  $\llbracket\mathrm{PT}_{\boldsymbol{\mathcal{A}}^{(2)}}\rrbracket_{s}$. If its $0$-composition is defined, i.e., if 
\[\mathrm{sc}^{0\llbracket\mathbf{PT}_{\boldsymbol{\mathcal{A}}^{(2)}}\rrbracket}_{s}\left(
\llbracket Q\rrbracket_{s}
\right)=
\mathrm{tg}^{0\llbracket\mathbf{PT}_{\boldsymbol{\mathcal{A}}^{(2)}}\rrbracket}_{s}\left(
\llbracket P\rrbracket_{s}
\right),
\]
then the following equalities hold
\allowdisplaybreaks
\begin{align*}
\mathrm{sc}^{0\llbracket\mathbf{PT}_{\boldsymbol{\mathcal{A}}^{(2)}}\rrbracket}_{s}\left(
\llbracket Q\rrbracket_{s}
\circ^{0\llbracket\mathbf{PT}_{\boldsymbol{\mathcal{A}}^{(2)}}\rrbracket}_{s}
\llbracket P\rrbracket_{s}
\right)
&=
\mathrm{sc}^{0\llbracket\mathbf{PT}_{\boldsymbol{\mathcal{A}}^{(2)}}\rrbracket}_{s}\left(
\llbracket P\rrbracket_{s}
\right);
\\
\mathrm{tg}^{0\llbracket\mathbf{PT}_{\boldsymbol{\mathcal{A}}^{(2)}}\rrbracket}_{s}\left(
\llbracket Q\rrbracket_{s}
\circ^{0\llbracket\mathbf{PT}_{\boldsymbol{\mathcal{A}}^{(2)}}\rrbracket}_{s}
\llbracket P\rrbracket_{s}
\right)
&=
\mathrm{tg}^{0\llbracket\mathbf{PT}_{\boldsymbol{\mathcal{A}}^{(2)}}\rrbracket}_{s}\left(
\llbracket P\rrbracket_{s}
\right).
\end{align*}
\end{proposition}
\begin{proof}
We will only present the proof of the first equality. The other one is handled in a similar manner.

The following chain of equalities holds
\begin{flushleft}
$\mathrm{sc}^{0\llbracket\mathbf{PT}_{\boldsymbol{\mathcal{A}}^{(2)}}\rrbracket}_{s}\left(
\llbracket Q\rrbracket_{s}
\circ^{0\llbracket\mathbf{PT}_{\boldsymbol{\mathcal{A}}^{(2)}}\rrbracket}_{s}
\llbracket P\rrbracket_{s}
\right)$
\allowdisplaybreaks
\begin{align*}
\qquad&=
\biggl\llbracket 
\mathrm{sc}^{0\mathbf{PT}_{\boldsymbol{\mathcal{A}}^{(2)}}}_{s}\left(
Q
\circ^{0\mathbf{PT}_{\boldsymbol{\mathcal{A}}^{(2)}}}_{s}
P
\right)
\biggr\rrbracket_{s}
\tag{1}
\\&=
\biggl\llbracket 
\mathrm{CH}^{(2)}_{s}\left(
\mathrm{ip}^{(2,X)@}_{s}\left(
\mathrm{sc}^{0\mathbf{PT}_{\boldsymbol{\mathcal{A}}^{(2)}}}_{s}\left(
Q
\circ^{0\mathbf{PT}_{\boldsymbol{\mathcal{A}}^{(2)}}}_{s}
P
\right)\right)\right)
\biggr\rrbracket_{s}
\tag{2}
\\&=
\biggl\llbracket 
\mathrm{CH}^{(2)}_{s}\left(
\mathrm{sc}^{0\mathbf{Pth}_{\boldsymbol{\mathcal{A}}^{(2)}}}_{s}\left(
\mathrm{ip}^{(2,X)@}_{s}\left(
Q
\right)
\circ^{0\mathbf{Pth}_{\boldsymbol{\mathcal{A}}^{(2)}}}_{s}
\mathrm{ip}^{(2,X)@}_{s}\left(
P
\right)
\right)\right)
\biggr\rrbracket_{s}
\tag{3}
\\&=
\biggl\llbracket 
\mathrm{CH}^{(2)}_{s}\left(
\mathrm{sc}^{0\mathbf{Pth}_{\boldsymbol{\mathcal{A}}^{(2)}}}_{s}\left(
\mathrm{ip}^{(2,X)@}_{s}\left(
P
\right)
\right)\right)
\biggr\rrbracket_{s}
\tag{4}
\\&=
\biggl\llbracket 
\mathrm{CH}^{(2)}_{s}\left(
\mathrm{ip}^{(2,X)@}_{s}\left(
\mathrm{sc}^{0\mathbf{PT}_{\boldsymbol{\mathcal{A}}^{(2)}}}_{s}\left(
P
\right)\right)\right)
\biggr\rrbracket_{s}
\tag{5}
\\&=
\biggl\llbracket 
\mathrm{sc}^{0\mathbf{PT}_{\boldsymbol{\mathcal{A}}^{(2)}}}_{s}\left(
P
\right)
\biggr\rrbracket_{s}
\tag{6}
\\&=
\mathrm{sc}^{0
\llbracket\mathbf{PT}_{\boldsymbol{\mathcal{A}}^{(2)}}\rrbracket}_{s}\left(
\bigl\llbracket 
P
\bigr\rrbracket_{s}
\right).
\tag{7}
\end{align*}
\end{flushleft}

In the just stated chain of equalities, the first equality unravels the description of the operation of $0$-source and $0$-composition in the many sorted partial $\Sigma^{\boldsymbol{\mathcal{A}}^{(2)}}$-algebra $\llbracket\mathbf{PT}_{\boldsymbol{\mathcal{A}}^{(2)}}\rrbracket$ according to Proposition~\ref{PDPTQDCatAlg}; the second equality follows from Proposition~\ref{LDVWCong}; the third equality follows from the fact that $\mathrm{ip}^{(2,X)@}$ is a $\Sigma^{\boldsymbol{\mathcal{A}}^{(2)}}$-homomorphism according to Definition~\ref{DDIp}; the fourth equality follows from Proposition~\ref{PDVDCH} and Proposition~\ref{PDVVarA4}; the fifth equality follows from the fact that $\mathrm{ip}^{(2,X)@}$ is a $\Sigma^{\boldsymbol{\mathcal{A}}^{(2)}}$-homomorphism according to Definition~\ref{DDIp}; the sixth equality follows from Proposition~\ref{LDVWCong}; finally, the last equality recovers the description of the operation of $0$-source in the many sorted partial $\Sigma^{\boldsymbol{\mathcal{A}}^{(2)}}$-algebra $\llbracket\mathbf{PT}_{\boldsymbol{\mathcal{A}}^{(2)}}\rrbracket$ according to Proposition~\ref{PDPTQDCatAlg}.

This completes the proof.
\end{proof}

\begin{proposition}\label{PDPTQVarA5} Let $s$ be a sort in $S$ and $\llbracket P\rrbracket_{s}$ a second-order path term class in  $\llbracket\mathrm{PT}_{\boldsymbol{\mathcal{A}}^{(2)}}\rrbracket_{s}$, then the following equalities hold
\allowdisplaybreaks
\begin{align*}
\llbracket P\rrbracket_{s}
\circ^{0\llbracket\mathbf{PT}_{\boldsymbol{\mathcal{A}}^{(2)}}\rrbracket}_{s}
\mathrm{sc}^{0\llbracket\mathbf{PT}_{\boldsymbol{\mathcal{A}}^{(2)}}\rrbracket}_{s}\left(
\llbracket P\rrbracket_{s}
\right)
&=
\llbracket P\rrbracket_{s};
\\
\mathrm{tg}^{0\llbracket\mathbf{PT}_{\boldsymbol{\mathcal{A}}^{(2)}}\rrbracket}_{s}\left(
\llbracket P\rrbracket_{s}
\right)
\circ^{0\llbracket\mathbf{PT}_{\boldsymbol{\mathcal{A}}^{(2)}}\rrbracket}_{s}
\llbracket P\rrbracket_{s}
&=
\llbracket P\rrbracket_{s};
\end{align*}
\end{proposition}
\begin{proof}
We will only present the proof for the first equality. The other one is handled in a similar manner.

The following chain of equalities holds
\begin{flushleft}
$\llbracket P\rrbracket_{s}
\circ^{0\llbracket\mathbf{PT}_{\boldsymbol{\mathcal{A}}^{(2)}}\rrbracket}_{s}
\mathrm{sc}^{0\llbracket\mathbf{PT}_{\boldsymbol{\mathcal{A}}^{(2)}}\rrbracket}_{s}\left(
\llbracket P\rrbracket_{s}
\right)$
\allowdisplaybreaks
\begin{align*}
\qquad&=
\biggl\llbracket 
P
\circ^{0\mathbf{PT}_{\boldsymbol{\mathcal{A}}^{(2)}}}_{s}
\mathrm{sc}^{0\mathbf{PT}_{\boldsymbol{\mathcal{A}}^{(2)}}}_{s}\left(
P
\right)
\biggr\rrbracket_{s}
\tag{1}
\\&=
\biggl\llbracket 
\mathrm{CH}^{(2)}_{s}\left(
\mathrm{ip}^{(2,X)@}_{s}\left(
P
\circ^{0\mathbf{PT}_{\boldsymbol{\mathcal{A}}^{(2)}}}_{s}
\mathrm{sc}^{0\mathbf{PT}_{\boldsymbol{\mathcal{A}}^{(2)}}}_{s}\left(
P
\right)
\right)\right)
\biggr\rrbracket_{s}
\tag{2}
\\&=
\biggl\llbracket 
\mathrm{CH}^{(2)}_{s}\left(
\mathrm{ip}^{(2,X)@}_{s}\left(
P
\right)
\circ^{0\mathbf{Pth}_{\boldsymbol{\mathcal{A}}^{(2)}}}_{s}
\mathrm{sc}^{0\mathbf{Pth}_{\boldsymbol{\mathcal{A}}^{(2)}}}_{s}\left(
\mathrm{ip}^{(2,X)@}_{s}\left(
P
\right)
\right)\right)
\biggr\rrbracket_{s}
\tag{3}
\\&=
\biggl\llbracket 
\mathrm{CH}^{(2)}_{s}\left(
\mathrm{ip}^{(2,X)@}_{s}\left(
P
\right)\right)
\biggr\rrbracket_{s}
\tag{4}
\\&=
\llbracket 
P
\rrbracket_{s}
.
\tag{5}
\end{align*}
\end{flushleft}

In the just stated chain of equalities, the first equality unravels the description of the operation of $0$-source and $0$-composition in the many sorted partial $\Sigma^{\boldsymbol{\mathcal{A}}^{(2)}}$-algebra $\llbracket\mathbf{PT}_{\boldsymbol{\mathcal{A}}^{(2)}}\rrbracket$ according to Proposition~\ref{PDPTQDCatAlg}; the second equality follows from Proposition~\ref{LDVWCong}; the third equality follows from the fact that $\mathrm{ip}^{(2,X)@}$ is a $\Sigma^{\boldsymbol{\mathcal{A}}^{(2)}}$-homomorphism according to Definition~\ref{DDIp}; the fourth equality follows from Proposition~\ref{PDVDCH} and Proposition~\ref{PDVVarA5}; the fifth equality follows from Proposition~\ref{LDVWCong}.

This completes the proof.
\end{proof}

\begin{proposition}\label{PDPTQVarA6} Let $s$ be a sort in $S$ and $\llbracket P\rrbracket_{s}$, $\llbracket Q\rrbracket_{s}$, $\llbracket R\rrbracket_{s}$ second-order path term classes in  $\llbracket\mathrm{PT}_{\boldsymbol{\mathcal{A}}^{(2)}}\rrbracket_{s}$ satisfying that
\begin{align*}
\mathrm{sc}^{0\llbracket\mathbf{PT}_{\boldsymbol{\mathcal{A}}^{(2)}}\rrbracket}_{s}\left(
\llbracket R\rrbracket_{s}
\right)&=
\mathrm{tg}^{0\llbracket\mathbf{PT}_{\boldsymbol{\mathcal{A}}^{(2)}}\rrbracket}_{s}\left(
\llbracket Q\rrbracket_{s}
\right);
\\
\mathrm{sc}^{0\llbracket\mathbf{PT}_{\boldsymbol{\mathcal{A}}^{(2)}}\rrbracket}_{s}\left(
\llbracket Q\rrbracket_{s}
\right)&=
\mathrm{tg}^{0\llbracket\mathbf{PT}_{\boldsymbol{\mathcal{A}}^{(2)}}\rrbracket}_{s}\left(
\llbracket P\rrbracket_{s}
\right),
\end{align*}
then the following equalities hold
\[
\llbracket R\rrbracket_{s}
\circ^{0\llbracket\mathbf{PT}_{\boldsymbol{\mathcal{A}}^{(2)}}\rrbracket}_{s}
\left(
\llbracket Q\rrbracket_{s}
\circ^{0\llbracket\mathbf{PT}_{\boldsymbol{\mathcal{A}}^{(2)}}\rrbracket}_{s}
\llbracket P\rrbracket_{s}
\right)
=
\left(
\llbracket R\rrbracket_{s}
\circ^{0\llbracket\mathbf{PT}_{\boldsymbol{\mathcal{A}}^{(2)}}\rrbracket}_{s}
\llbracket Q\rrbracket_{s}
\right)
\circ^{0\llbracket\mathbf{PT}_{\boldsymbol{\mathcal{A}}^{(2)}}\rrbracket}_{s}
\llbracket P\rrbracket_{s}.
\]
\end{proposition}
\begin{proof}
The following chain of equalities holds
\begin{flushleft}
$\llbracket R\rrbracket_{s}
\circ^{0\llbracket\mathbf{PT}_{\boldsymbol{\mathcal{A}}^{(2)}}\rrbracket}_{s}
\left(
\llbracket Q\rrbracket_{s}
\circ^{0\llbracket\mathbf{PT}_{\boldsymbol{\mathcal{A}}^{(2)}}\rrbracket}_{s}
\llbracket P\rrbracket_{s}
\right)$
\allowdisplaybreaks
\begin{align*}
&=
\biggl\llbracket 
R
\circ^{0\mathbf{PT}_{\boldsymbol{\mathcal{A}}^{(2)}}}_{s}
\left(
Q
\circ^{0\mathbf{PT}_{\boldsymbol{\mathcal{A}}^{(2)}}}_{s}
P
\right)
\biggr\rrbracket_{s}
\tag{1}
\\&=
\biggl\llbracket 
\mathrm{CH}^{(2)}_{s}\left(
\mathrm{ip}^{(2,X)@}_{s}\left(
R
\circ^{0\mathbf{PT}_{\boldsymbol{\mathcal{A}}^{(2)}}}_{s}
\left(
Q
\circ^{0\mathbf{PT}_{\boldsymbol{\mathcal{A}}^{(2)}}}_{s}
P
\right)
\right)\right)
\biggr\rrbracket_{s}
\tag{2}
\\&=
\biggl\llbracket 
\mathrm{CH}^{(2)}_{s}\left(
\mathrm{ip}^{(2,X)@}_{s}\left(
R
\right)
\circ^{0\mathbf{Pth}_{\boldsymbol{\mathcal{A}}^{(2)}}}_{s}
\left(
\mathrm{ip}^{(2,X)@}_{s}\left(
Q
\right)
\circ^{0\mathbf{Pth}_{\boldsymbol{\mathcal{A}}^{(2)}}}_{s}
\mathrm{ip}^{(2,X)@}_{s}\left(
P
\right)
\right)
\right)
\biggr\rrbracket_{s}
\tag{3}
\\&=
\biggl\llbracket 
\mathrm{CH}^{(2)}_{s}\left(
\left(
\mathrm{ip}^{(2,X)@}_{s}\left(
R
\right)
\circ^{0\mathbf{Pth}_{\boldsymbol{\mathcal{A}}^{(2)}}}_{s}
\mathrm{ip}^{(2,X)@}_{s}\left(
Q
\right)
\right)
\circ^{0\mathbf{Pth}_{\boldsymbol{\mathcal{A}}^{(2)}}}_{s}
\mathrm{ip}^{(2,X)@}_{s}\left(
P
\right)
\right)
\biggr\rrbracket_{s}
\tag{4}
\\&=
\biggl\llbracket 
\mathrm{CH}^{(2)}_{s}\left(
\mathrm{ip}^{(2,X)@}_{s}\left(
\left(
R
\circ^{0\mathbf{PT}_{\boldsymbol{\mathcal{A}}^{(2)}}}_{s}
Q
\right)
\circ^{0\mathbf{PT}_{\boldsymbol{\mathcal{A}}^{(2)}}}_{s}
P
\right)\right)
\biggr\rrbracket_{s}
\tag{5}
\\&=
\biggl\llbracket 
\left(
R
\circ^{0\mathbf{PT}_{\boldsymbol{\mathcal{A}}^{(2)}}}_{s}
Q
\right)
\circ^{0\mathbf{PT}_{\boldsymbol{\mathcal{A}}^{(2)}}}_{s}
P
\biggr\rrbracket_{s}
\tag{6}
\\&=
\left(
\llbracket R\rrbracket_{s}
\circ^{0\llbracket\mathbf{PT}_{\boldsymbol{\mathcal{A}}^{(2)}}\rrbracket}_{s}
\llbracket Q\rrbracket_{s}
\right)
\circ^{0\llbracket\mathbf{PT}_{\boldsymbol{\mathcal{A}}^{(2)}}\rrbracket}_{s}
\llbracket P\rrbracket_{s}.
\tag{7}
\end{align*}
\end{flushleft}

In the just stated chain of equalities, the first equality unravels the description of the operation of $0$-composition in the many sorted partial $\Sigma^{\boldsymbol{\mathcal{A}}^{(2)}}$-algebra $\llbracket\mathbf{PT}_{\boldsymbol{\mathcal{A}}^{(2)}}\rrbracket$ according to Proposition~\ref{PDPTQDCatAlg}; the second equality follows from Proposition~\ref{LDVWCong}; the third equality follows from the fact that $\mathrm{ip}^{(2,X)@}$ is a $\Sigma^{\boldsymbol{\mathcal{A}}^{(2)}}$-homomorphism according to Definition~\ref{DDIp}; the fourth equality follows from Proposition~\ref{PDVDCH} and Proposition~\ref{PDVVarA6}; the fifth equality follows from the fact that $\mathrm{ip}^{(2,X)@}$ is a $\Sigma^{\boldsymbol{\mathcal{A}}^{(2)}}$-homomorphism according to Definition~\ref{DDIp}; the sixth equality follows from Proposition~\ref{LDVWCong}; finally, the last equality recovers the description of the operation of $0$-composition in the many sorted partial $\Sigma^{\boldsymbol{\mathcal{A}}^{(2)}}$-algebra $\llbracket\mathbf{PT}_{\boldsymbol{\mathcal{A}}^{(2)}}\rrbracket$ according to Proposition~\ref{PDPTQDCatAlg}.

This completes the proof.
\end{proof}

\begin{proposition}\label{PDPTQVarB2} Let $s$ be a sort in $S$ and $\llbracket P\rrbracket_{s}$ a second-order path term class in  $\llbracket\mathrm{PT}_{\boldsymbol{\mathcal{A}}^{(2)}}\rrbracket_{s}$, then the following equalities hold
\allowdisplaybreaks
\begin{align*}
\mathrm{sc}^{1\llbracket\mathbf{PT}_{\boldsymbol{\mathcal{A}}^{(2)}}\rrbracket}_{s}\left(
\mathrm{sc}^{1\llbracket\mathbf{PT}_{\boldsymbol{\mathcal{A}}^{(2)}}\rrbracket}_{s}\left(
\llbracket P\rrbracket_{s}
\right)
\right)
& 
=
\mathrm{sc}^{1\llbracket\mathbf{PT}_{\boldsymbol{\mathcal{A}}^{(2)}}\rrbracket}_{s}\left(
\llbracket P\rrbracket_{s}
\right);
\\
\mathrm{sc}^{1\llbracket\mathbf{PT}_{\boldsymbol{\mathcal{A}}^{(2)}}\rrbracket}_{s}\left(
\mathrm{tg}^{1\llbracket\mathbf{PT}_{\boldsymbol{\mathcal{A}}^{(2)}}\rrbracket}_{s}\left(
\llbracket P\rrbracket_{s}
\right)
\right)
& 
=
\mathrm{tg}^{1\llbracket\mathbf{PT}_{\boldsymbol{\mathcal{A}}^{(2)}}\rrbracket}_{s}\left(
\llbracket P\rrbracket_{s}
\right);
\\
\mathrm{tg}^{1\llbracket\mathbf{PT}_{\boldsymbol{\mathcal{A}}^{(2)}}\rrbracket}_{s}\left(
\mathrm{sc}^{1\llbracket\mathbf{PT}_{\boldsymbol{\mathcal{A}}^{(2)}}\rrbracket}_{s}\left(
\llbracket P\rrbracket_{s}
\right)
\right)
& 
=
\mathrm{sc}^{1\llbracket\mathbf{PT}_{\boldsymbol{\mathcal{A}}^{(2)}}\rrbracket}_{s}\left(
\llbracket P\rrbracket_{s}
\right);
\\
\mathrm{tg}^{1\llbracket\mathbf{PT}_{\boldsymbol{\mathcal{A}}^{(2)}}\rrbracket}_{s}\left(
\mathrm{tg}^{1\llbracket\mathbf{PT}_{\boldsymbol{\mathcal{A}}^{(2)}}\rrbracket}_{s}\left(
\llbracket P\rrbracket_{s}
\right)
\right)
& 
=
\mathrm{tg}^{1\llbracket\mathbf{PT}_{\boldsymbol{\mathcal{A}}^{(2)}}\rrbracket}_{s}\left(
\llbracket P\rrbracket_{s}
\right).
\end{align*}	
\end{proposition}
\begin{proof}
We will only present the proof for the first equality. The other three are handled in a similar manner.

Note that the following chain of equalities holds
\begin{flushleft}
$\mathrm{sc}^{1\llbracket\mathbf{PT}_{\boldsymbol{\mathcal{A}}^{(2)}}\rrbracket}_{s}\left(
\mathrm{sc}^{1\llbracket\mathbf{PT}_{\boldsymbol{\mathcal{A}}^{(2)}}\rrbracket}_{s}\left(
\llbracket P\rrbracket_{s}
\right)
\right)$
\allowdisplaybreaks
\begin{align*}
\quad&=
\biggl\llbracket
\mathrm{sc}^{1\mathbf{PT}_{\boldsymbol{\mathcal{A}}^{(2)}}}_{s}\left(
\mathrm{sc}^{1\mathbf{PT}_{\boldsymbol{\mathcal{A}}^{(2)}}}_{s}\left(
P
\right)\right)
\biggr\rrbracket_{s}
\tag{1}
\\&=
\biggl\llbracket
\mathrm{CH}^{(2)}_{s}\left(
\mathrm{ip}^{(2,X)@}_{s}\left(
\mathrm{sc}^{1\mathbf{PT}_{\boldsymbol{\mathcal{A}}^{(2)}}}_{s}\left(
\mathrm{sc}^{1\mathbf{PT}_{\boldsymbol{\mathcal{A}}^{(2)}}}_{s}\left(
P
\right)\right)
\right)\right)
\biggr\rrbracket_{s}
\tag{2}
\\&=
\biggl\llbracket
\mathrm{CH}^{(2)}_{s}\left(
\mathrm{sc}^{1\mathbf{Pth}_{\boldsymbol{\mathcal{A}}^{(2)}}}_{s}\left(
\mathrm{sc}^{1\mathbf{Pth}_{\boldsymbol{\mathcal{A}}^{(2)}}}_{s}\left(
\mathrm{ip}^{(2,X)@}_{s}\left(
P
\right)\right)
\right)\right)
\biggr\rrbracket_{s}
\tag{3}
\\&=
\biggl\llbracket
\mathrm{CH}^{(2)}_{s}\left(
\mathrm{sc}^{1\mathbf{Pth}_{\boldsymbol{\mathcal{A}}^{(2)}}}_{s}\left(
\mathrm{ip}^{(2,X)@}_{s}\left(
P
\right)
\right)\right)
\biggr\rrbracket_{s}
\tag{4}
\\&=
\biggl\llbracket
\mathrm{CH}^{(2)}_{s}\left(
\mathrm{ip}^{(2,X)@}_{s}\left(
\mathrm{sc}^{1\mathbf{PT}_{\boldsymbol{\mathcal{A}}^{(2)}}}_{s}\left(
P
\right)
\right)\right)
\biggr\rrbracket_{s}
\tag{5}
\\&=
\biggl\llbracket
\mathrm{sc}^{1\mathbf{PT}_{\boldsymbol{\mathcal{A}}^{(2)}}}_{s}\left(
P
\right)
\biggr\rrbracket_{s}
\tag{6}
\\&=
\mathrm{sc}^{1\llbracket\mathbf{PT}_{\boldsymbol{\mathcal{A}}^{(2)}}\rrbracket}_{s}\left(
\llbracket P\rrbracket_{s}
\right).
\tag{7}
\end{align*}
\end{flushleft}

In the just stated chain of equalities, the first equality unravels the description of the operation of $1$-source in the many sorted partial $\Sigma^{\boldsymbol{\mathcal{A}}^{(2)}}$-algebra $\llbracket\mathbf{PT}_{\boldsymbol{\mathcal{A}}^{(2)}}\rrbracket$ according to Proposition~\ref{PDPTQDCatAlg}; the second equality follows from Proposition~\ref{LDVWCong}; the third equality follows from the fact that $\mathrm{ip}^{(2,X)@}$ is a $\Sigma^{\boldsymbol{\mathcal{A}}^{(2)}}$-homomorphism according to Definition~\ref{DDIp}; the fourth equality follows from Proposition~\ref{PDVDCH} and Proposition~\ref{PDVVarB2}; the fifth equality follows from the fact that $\mathrm{ip}^{(2,X)@}$ is a $\Sigma^{\boldsymbol{\mathcal{A}}^{(2)}}$-homomorphism according to Definition~\ref{DDIp}; the sixth equality follows from Proposition~\ref{LDVWCong}; finally, the last equality recovers the description of the operation of $1$-source in the many sorted partial $\Sigma^{\boldsymbol{\mathcal{A}}^{(2)}}$-algebra $\llbracket\mathbf{PT}_{\boldsymbol{\mathcal{A}}^{(2)}}\rrbracket$ according to Proposition~\ref{PDPTQDCatAlg}.

This completes the proof.
\end{proof}

\begin{proposition}\label{PDPTQVarB3} Let $s$ be a sort in $S$ and $\llbracket P\rrbracket_{s}$, $\llbracket Q\rrbracket_{s}$ second-order path term classes in  $\llbracket\mathrm{PT}_{\boldsymbol{\mathcal{A}}^{(2)}}\rrbracket_{s}$, then the following statement are equivalent
\begin{itemize}
\item[(i)] $\llbracket Q\rrbracket_{s}
\circ^{1\llbracket\mathbf{PT}_{\boldsymbol{\mathcal{A}}^{(2)}}\rrbracket}_{s}
\llbracket P\rrbracket_{s}
$ is defined;
\item[(ii)] $\mathrm{sc}^{1\llbracket\mathbf{PT}_{\boldsymbol{\mathcal{A}}^{(2)}}\rrbracket}_{s}(
\llbracket Q\rrbracket_{s}
)=
\mathrm{tg}^{1\llbracket\mathbf{PT}_{\boldsymbol{\mathcal{A}}^{(2)}}\rrbracket}_{s}(
\llbracket P\rrbracket_{s}
)
$.
\end{itemize}
\end{proposition}
\begin{proof}
The following chain of equivalences holds
\begin{flushleft}
$\llbracket Q\rrbracket_{s}
\circ^{1\llbracket\mathbf{PT}_{\boldsymbol{\mathcal{A}}^{(2)}}\rrbracket}_{s}
\llbracket P\rrbracket_{s}$
\mbox{ is defined}
\allowdisplaybreaks
\begin{align*}
&\Leftrightarrow
\mathrm{sc}^{([1],2)}_{s}\left(\mathrm{ip}^{(2,X)@}_{s}\left(
Q\right)\right)
=
\mathrm{tg}^{([1],2)}_{s}\left(\mathrm{ip}^{(2,X)@}_{s}\left(
P\right)\right)
\tag{1}
\\&\Leftrightarrow
\mathrm{sc}^{1\llbracket \mathbf{Pth}_{\boldsymbol{\mathcal{A}}^{(2)}} \rrbracket}_{s}\left(
\Bigl\llbracket \mathrm{ip}^{(2,X)@}_{s}\left(Q\right) \Bigr\rrbracket_{s}\right)
=
\mathrm{tg}^{1\llbracket \mathbf{Pth}_{\boldsymbol{\mathcal{A}}^{(2)}} \rrbracket}_{s}\left(
\Bigl\llbracket \mathrm{ip}^{(2,X)@}_{s}\left(P\right) \Bigr\rrbracket_{s}\right)
\tag{2}
\\&\Leftrightarrow
\Bigl\llbracket \mathrm{sc}^{1\mathbf{Pth}_{\boldsymbol{\mathcal{A}}^{(2)}}}_{s}\left(\mathrm{ip}^{(2,X)@}_{s}\left(Q\right)\right) \Bigr\rrbracket_{s}
=
\Bigl\llbracket \mathrm{tg}^{1\mathbf{Pth}_{\boldsymbol{\mathcal{A}}^{(2)}}}_{s}\left(\mathrm{ip}^{(2,X)@}_{s}\left(Q\right)\right) \Bigr\rrbracket_{s}
\tag{3}
\\&\Leftrightarrow
\biggl\llbracket \mathrm{ip}^{(2,X)@}_{s}\left( \mathrm{sc}^{1\mathbf{PT}_{\boldsymbol{\mathcal{A}}^{(2)}}}_{s}\left(Q\right)\right) \biggr\rrbracket_{s}=
\biggl\llbracket \mathrm{ip}^{(2,X)@}_{s}\left( \mathrm{tg}^{1\mathbf{PT}_{\boldsymbol{\mathcal{A}}^{(2)}}}_{s}\left(Q\right)\right) \biggr\rrbracket_{s}
\tag{4}
\\&\Leftrightarrow
\biggl\llbracket \mathrm{CH}^{(2)}_{s}\left(\mathrm{ip}^{(2,X)@}_{s}\left( \mathrm{sc}^{1\mathbf{PT}_{\boldsymbol{\mathcal{A}}^{(2)}}}_{s}\left(Q\right)\right)\right) \biggr\rrbracket_{s}
\\&\qquad\qquad\qquad\qquad=
\biggl\llbracket \mathrm{CH}^{(2)}_{s}\left(\mathrm{ip}^{(2,X)@}_{s}\left( \mathrm{tg}^{1\mathbf{PT}_{\boldsymbol{\mathcal{A}}^{(2)}}}_{s}\left(Q\right)\right)\right) \biggr\rrbracket_{s}
\tag{5}
\\&\Leftrightarrow
\Bigl\llbracket  \mathrm{sc}^{1\mathbf{PT}_{\boldsymbol{\mathcal{A}}^{(2)}}}_{s}\left(Q\right) \Bigr\rrbracket_{s}
=
\Bigl\llbracket  \mathrm{tg}^{1\mathbf{PT}_{\boldsymbol{\mathcal{A}}^{(2)}}}_{s}\left(P\right) \Bigr\rrbracket_{s}
\tag{6}
\\&\Leftrightarrow
\mathrm{sc}^{1\llbracket\mathbf{PT}_{\boldsymbol{\mathcal{A}}^{(2)}}\rrbracket}_{s}\left(
\llbracket Q\rrbracket_{s}
\right)=
\mathrm{tg}^{1\llbracket\mathbf{PT}_{\boldsymbol{\mathcal{A}}^{(2)}}\rrbracket}_{s}\left(
\llbracket P\rrbracket_{s}
\right).
\tag{7}
\end{align*}
\end{flushleft}

In the just stated chain of equivalences, the first equivalence follows from the description of the $1$-composition operation in the many-sorted partial $\Sigma^{\boldsymbol{\mathcal{A}}^{(2)}}$-algebra $\llbracket\mathbf{PT}_{\boldsymbol{\mathcal{A}}^{(2)}}\rrbracket$ according to Proposition~\ref{PDPTQDCatAlg}; the second equivalence follows from Proposition~\ref{PDVVarB3}; the third equivalence simply unravels the description of the $1$-source and $1$-target operations in the many-sorted partial $\Sigma^{\boldsymbol{\mathcal{A}}^{(2)}}$-algebra $\llbracket \mathbf{Pth}_{\boldsymbol{\mathcal{A}}^{(2)}} \rrbracket$, according to Proposition~\ref{PDVDCatAlg}; the fourth equivalence follows from the fact that, according to Definition~\ref{DDIp}, $\mathrm{ip}^{(2,X)@}$ is a $\Sigma^{\boldsymbol{\mathcal{A}}^{(2)}}$-homomorphism; the fifth equivalence follows from left to right from Proposition~\ref{PDVDCH} and from right to left by Corollary~\ref{CDVCongDCH}; the sixth equivalence follows from Lemma~\ref{LDVWCong}; finally, the last equivalence recovers the description of the $1$-source and $1$-target operations in the many-sorted partial $\Sigma^{\boldsymbol{\mathcal{A}}^{(2)}}$-algebra $\llbracket\mathbf{PT}_{\boldsymbol{\mathcal{A}}^{(2)}}\rrbracket$ according to Proposition~\ref{PDPTQDCatAlg}.

This completes the proof.
\end{proof}

\begin{proposition}\label{PDPTQVarB4} Let $s$ be a sort in $S$ and $\llbracket P\rrbracket_{s}$, $\llbracket Q\rrbracket_{s}$ second-order path term classes in  $\llbracket\mathrm{PT}_{\boldsymbol{\mathcal{A}}^{(2)}}\rrbracket_{s}$. If its $1$-composition is defined, i.e., if 
\[\mathrm{sc}^{1\llbracket\mathbf{PT}_{\boldsymbol{\mathcal{A}}^{(2)}}\rrbracket}_{s}\left(
\llbracket Q\rrbracket_{s}
\right)=
\mathrm{tg}^{1\llbracket\mathbf{PT}_{\boldsymbol{\mathcal{A}}^{(2)}}\rrbracket}_{s}\left(
\llbracket P\rrbracket_{s}
\right),
\]
then the following equalities hold
\allowdisplaybreaks
\begin{align*}
\mathrm{sc}^{1\llbracket\mathbf{PT}_{\boldsymbol{\mathcal{A}}^{(2)}}\rrbracket}_{s}\left(
\llbracket Q\rrbracket_{s}
\circ^{1\llbracket\mathbf{PT}_{\boldsymbol{\mathcal{A}}^{(2)}}\rrbracket}_{s}
\llbracket P\rrbracket_{s}
\right)
&=
\mathrm{sc}^{1\llbracket\mathbf{PT}_{\boldsymbol{\mathcal{A}}^{(2)}}\rrbracket}_{s}\left(
\llbracket P\rrbracket_{s}
\right);
\\
\mathrm{tg}^{1\llbracket\mathbf{PT}_{\boldsymbol{\mathcal{A}}^{(2)}}\rrbracket}_{s}\left(
\llbracket Q\rrbracket_{s}
\circ^{1\llbracket\mathbf{PT}_{\boldsymbol{\mathcal{A}}^{(2)}}\rrbracket}_{s}
\llbracket P\rrbracket_{s}
\right)
&=
\mathrm{tg}^{1\llbracket\mathbf{PT}_{\boldsymbol{\mathcal{A}}^{(2)}}\rrbracket}_{s}\left(
\llbracket P\rrbracket_{s}
\right).
\end{align*}
\end{proposition}
\begin{proof}
We will only present the proof of the first equality. The other one is handled in a similar manner.

The following chain of equalities holds
\begin{flushleft}
$\mathrm{sc}^{1\llbracket\mathbf{PT}_{\boldsymbol{\mathcal{A}}^{(2)}}\rrbracket}_{s}\left(
\llbracket Q\rrbracket_{s}
\circ^{1\llbracket\mathbf{PT}_{\boldsymbol{\mathcal{A}}^{(2)}}\rrbracket}_{s}
\llbracket P\rrbracket_{s}
\right)$
\allowdisplaybreaks
\begin{align*}
\qquad&=
\biggl\llbracket 
\mathrm{sc}^{1\mathbf{PT}_{\boldsymbol{\mathcal{A}}^{(2)}}}_{s}\left(
Q
\circ^{1\mathbf{PT}_{\boldsymbol{\mathcal{A}}^{(2)}}}_{s}
P
\right)
\biggr\rrbracket_{s}
\tag{1}
\\&=
\biggl\llbracket 
\mathrm{CH}^{(2)}_{s}\left(
\mathrm{ip}^{(2,X)@}_{s}\left(
\mathrm{sc}^{1\mathbf{PT}_{\boldsymbol{\mathcal{A}}^{(2)}}}_{s}\left(
Q
\circ^{1\mathbf{PT}_{\boldsymbol{\mathcal{A}}^{(2)}}}_{s}
P
\right)\right)\right)
\biggr\rrbracket_{s}
\tag{2}
\\&=
\biggl\llbracket 
\mathrm{CH}^{(2)}_{s}\left(
\mathrm{sc}^{1\mathbf{Pth}_{\boldsymbol{\mathcal{A}}^{(2)}}}_{s}\left(
\mathrm{ip}^{(2,X)@}_{s}\left(
Q
\right)
\circ^{1\mathbf{Pth}_{\boldsymbol{\mathcal{A}}^{(2)}}}_{s}
\mathrm{ip}^{(2,X)@}_{s}\left(
P
\right)
\right)\right)
\biggr\rrbracket_{s}
\tag{3}
\\&=
\biggl\llbracket 
\mathrm{CH}^{(2)}_{s}\left(
\mathrm{sc}^{1\mathbf{Pth}_{\boldsymbol{\mathcal{A}}^{(2)}}}_{s}\left(
\mathrm{ip}^{(2,X)@}_{s}\left(
P
\right)
\right)\right)
\biggr\rrbracket_{s}
\tag{4}
\\&=
\biggl\llbracket 
\mathrm{CH}^{(2)}_{s}\left(
\mathrm{ip}^{(2,X)@}_{s}\left(
\mathrm{sc}^{1\mathbf{PT}_{\boldsymbol{\mathcal{A}}^{(2)}}}_{s}\left(
P
\right)\right)\right)
\biggr\rrbracket_{s}
\tag{5}
\\&=
\biggl\llbracket 
\mathrm{sc}^{1\mathbf{PT}_{\boldsymbol{\mathcal{A}}^{(2)}}}_{s}\left(
P
\right)
\biggr\rrbracket_{s}
\tag{6}
\\&=
\mathrm{sc}^{1
\llbracket\mathbf{PT}_{\boldsymbol{\mathcal{A}}^{(2)}}\rrbracket}_{s}\left(
\bigl\llbracket 
P
\bigr\rrbracket_{s}
\right).
\tag{7}
\end{align*}
\end{flushleft}

In the just stated chain of equalities, the first equality unravels the description of the operation of $1$-source and $1$-composition in the many sorted partial $\Sigma^{\boldsymbol{\mathcal{A}}^{(2)}}$-algebra $\llbracket\mathbf{PT}_{\boldsymbol{\mathcal{A}}^{(2)}}\rrbracket$ according to Proposition~\ref{PDPTQDCatAlg}; the second equality follows from Proposition~\ref{LDVWCong}; the third equality follows from the fact that $\mathrm{ip}^{(2,X)@}$ is a $\Sigma^{\boldsymbol{\mathcal{A}}^{(2)}}$-homomorphism according to Definition~\ref{DDIp}; the fourth equality follows from Proposition~\ref{PDVDCH} and Proposition~\ref{PDVVarB4}; the fifth equality follows from the fact that $\mathrm{ip}^{(2,X)@}$ is a $\Sigma^{\boldsymbol{\mathcal{A}}^{(2)}}$-homomorphism according to Definition~\ref{DDIp}; the sixth equality follows from Proposition~\ref{LDVWCong}; finally, the last equality recovers the description of the operation of $1$-source in the many sorted partial $\Sigma^{\boldsymbol{\mathcal{A}}^{(2)}}$-algebra $\llbracket\mathbf{PT}_{\boldsymbol{\mathcal{A}}^{(2)}}\rrbracket$ according to Proposition~\ref{PDPTQDCatAlg}.

This completes the proof.
\end{proof}

\begin{proposition}\label{PDPTQVarB5} Let $s$ be a sort in $S$ and $\llbracket P\rrbracket_{s}$ a second-order path term class in  $\llbracket\mathrm{PT}_{\boldsymbol{\mathcal{A}}^{(2)}}\rrbracket_{s}$, then the following equalities hold
\allowdisplaybreaks
\begin{align*}
\llbracket P\rrbracket_{s}
\circ^{1\llbracket\mathbf{PT}_{\boldsymbol{\mathcal{A}}^{(2)}}\rrbracket}_{s}
\mathrm{sc}^{1\llbracket\mathbf{PT}_{\boldsymbol{\mathcal{A}}^{(2)}}\rrbracket}_{s}\left(
\llbracket P\rrbracket_{s}
\right)
&=
\llbracket P\rrbracket_{s};
\\
\mathrm{tg}^{1\llbracket\mathbf{PT}_{\boldsymbol{\mathcal{A}}^{(2)}}\rrbracket}_{s}\left(
\llbracket P\rrbracket_{s}
\right)
\circ^{1\llbracket\mathbf{PT}_{\boldsymbol{\mathcal{A}}^{(2)}}\rrbracket}_{s}
\llbracket P\rrbracket_{s}
&=
\llbracket P\rrbracket_{s};
\end{align*}
\end{proposition}
\begin{proof}
We will only present the proof for the first equality. The other one is handled in a similar manner.

The following chain of equalities holds
\begin{flushleft}
$\llbracket P\rrbracket_{s}
\circ^{1\llbracket\mathbf{PT}_{\boldsymbol{\mathcal{A}}^{(2)}}\rrbracket}_{s}
\mathrm{sc}^{1\llbracket\mathbf{PT}_{\boldsymbol{\mathcal{A}}^{(2)}}\rrbracket}_{s}\left(
\llbracket P\rrbracket_{s}
\right)$
\allowdisplaybreaks
\begin{align*}
\qquad&=
\biggl\llbracket 
P
\circ^{1\mathbf{PT}_{\boldsymbol{\mathcal{A}}^{(2)}}}_{s}
\mathrm{sc}^{1\mathbf{PT}_{\boldsymbol{\mathcal{A}}^{(2)}}}_{s}\left(
P
\right)
\biggr\rrbracket_{s}
\tag{1}
\\&=
\biggl\llbracket 
\mathrm{CH}^{(2)}_{s}\left(
\mathrm{ip}^{(2,X)@}_{s}\left(
P
\circ^{1\mathbf{PT}_{\boldsymbol{\mathcal{A}}^{(2)}}}_{s}
\mathrm{sc}^{1\mathbf{PT}_{\boldsymbol{\mathcal{A}}^{(2)}}}_{s}\left(
P
\right)
\right)\right)
\biggr\rrbracket_{s}
\tag{2}
\\&=
\biggl\llbracket 
\mathrm{CH}^{(2)}_{s}\left(
\mathrm{ip}^{(2,X)@}_{s}\left(
P
\right)
\circ^{1\mathbf{Pth}_{\boldsymbol{\mathcal{A}}^{(2)}}}_{s}
\mathrm{sc}^{1\mathbf{Pth}_{\boldsymbol{\mathcal{A}}^{(2)}}}_{s}\left(
\mathrm{ip}^{(2,X)@}_{s}\left(
P
\right)
\right)\right)
\biggr\rrbracket_{s}
\tag{3}
\\&=
\biggl\llbracket 
\mathrm{CH}^{(2)}_{s}\left(
\mathrm{ip}^{(2,X)@}_{s}\left(
P
\right)\right)
\biggr\rrbracket_{s}
\tag{4}
\\&=
\llbracket 
P
\rrbracket_{s}
.
\tag{5}
\end{align*}
\end{flushleft}

In the just stated chain of equalities, the first equality unravels the description of the operation of $1$-source and $1$-composition in the many sorted partial $\Sigma^{\boldsymbol{\mathcal{A}}^{(2)}}$-algebra $\llbracket\mathbf{PT}_{\boldsymbol{\mathcal{A}}^{(2)}}\rrbracket$ according to Proposition~\ref{PDPTQDCatAlg}; the second equality follows from Proposition~\ref{LDVWCong}; the third equality follows from the fact that $\mathrm{ip}^{(2,X)@}$ is a $\Sigma^{\boldsymbol{\mathcal{A}}^{(2)}}$-homomorphism according to Definition~\ref{DDIp}; the fourth equality follows from Proposition~\ref{PDVDCH} and Proposition~\ref{PDVVarB5}; the fifth equality follows from Proposition~\ref{LDVWCong}.

This completes the proof.
\end{proof}

\begin{proposition}\label{PDPTQVarB6} Let $s$ be a sort in $S$ and $\llbracket P\rrbracket_{s}$, $\llbracket Q\rrbracket_{s}$, $\llbracket R\rrbracket_{s}$ second-order path term classes in  $\llbracket\mathrm{PT}_{\boldsymbol{\mathcal{A}}^{(2)}}\rrbracket_{s}$ satisfying that
\begin{align*}
\mathrm{sc}^{1\llbracket\mathbf{PT}_{\boldsymbol{\mathcal{A}}^{(2)}}\rrbracket}_{s}\left(
\llbracket R\rrbracket_{s}
\right)&=
\mathrm{tg}^{1\llbracket\mathbf{PT}_{\boldsymbol{\mathcal{A}}^{(2)}}\rrbracket}_{s}\left(
\llbracket Q\rrbracket_{s}
\right);
\\
\mathrm{sc}^{1\llbracket\mathbf{PT}_{\boldsymbol{\mathcal{A}}^{(2)}}\rrbracket}_{s}\left(
\llbracket Q\rrbracket_{s}
\right)&=
\mathrm{tg}^{1\llbracket\mathbf{PT}_{\boldsymbol{\mathcal{A}}^{(2)}}\rrbracket}_{s}\left(
\llbracket P\rrbracket_{s}
\right),
\end{align*}
then the following equalities hold
\[
\llbracket R\rrbracket_{s}
\circ^{1\llbracket\mathbf{PT}_{\boldsymbol{\mathcal{A}}^{(2)}}\rrbracket}_{s}
\left(
\llbracket Q\rrbracket_{s}
\circ^{1\llbracket\mathbf{PT}_{\boldsymbol{\mathcal{A}}^{(2)}}\rrbracket}_{s}
\llbracket P\rrbracket_{s}
\right)
=
\left(
\llbracket R\rrbracket_{s}
\circ^{1\llbracket\mathbf{PT}_{\boldsymbol{\mathcal{A}}^{(2)}}\rrbracket}_{s}
\llbracket Q\rrbracket_{s}
\right)
\circ^{1\llbracket\mathbf{PT}_{\boldsymbol{\mathcal{A}}^{(2)}}\rrbracket}_{s}
\llbracket P\rrbracket_{s}.
\]
\end{proposition}
\begin{proof}
The following chain of equalities holds
\begin{flushleft}
$\llbracket R\rrbracket_{s}
\circ^{1\llbracket\mathbf{PT}_{\boldsymbol{\mathcal{A}}^{(2)}}\rrbracket}_{s}
\left(
\llbracket Q\rrbracket_{s}
\circ^{1\llbracket\mathbf{PT}_{\boldsymbol{\mathcal{A}}^{(2)}}\rrbracket}_{s}
\llbracket P\rrbracket_{s}
\right)$
\allowdisplaybreaks
\begin{align*}
&=
\biggl\llbracket 
R
\circ^{1\mathbf{PT}_{\boldsymbol{\mathcal{A}}^{(2)}}}_{s}
\left(
Q
\circ^{1\mathbf{PT}_{\boldsymbol{\mathcal{A}}^{(2)}}}_{s}
P
\right)
\biggr\rrbracket_{s}
\tag{1}
\\&=
\biggl\llbracket 
\mathrm{CH}^{(2)}_{s}\left(
\mathrm{ip}^{(2,X)@}_{s}\left(
R
\circ^{1\mathbf{PT}_{\boldsymbol{\mathcal{A}}^{(2)}}}_{s}
\left(
Q
\circ^{1\mathbf{PT}_{\boldsymbol{\mathcal{A}}^{(2)}}}_{s}
P
\right)
\right)\right)
\biggr\rrbracket_{s}
\tag{2}
\\&=
\biggl\llbracket 
\mathrm{CH}^{(2)}_{s}\left(
\mathrm{ip}^{(2,X)@}_{s}\left(
R
\right)
\circ^{1\mathbf{Pth}_{\boldsymbol{\mathcal{A}}^{(2)}}}_{s}
\left(
\mathrm{ip}^{(2,X)@}_{s}\left(
Q
\right)
\circ^{1\mathbf{Pth}_{\boldsymbol{\mathcal{A}}^{(2)}}}_{s}
\mathrm{ip}^{(2,X)@}_{s}\left(
P
\right)
\right)
\right)
\biggr\rrbracket_{s}
\tag{3}
\\&=
\biggl\llbracket 
\mathrm{CH}^{(2)}_{s}\left(
\left(
\mathrm{ip}^{(2,X)@}_{s}\left(
R
\right)
\circ^{1\mathbf{Pth}_{\boldsymbol{\mathcal{A}}^{(2)}}}_{s}
\mathrm{ip}^{(2,X)@}_{s}\left(
Q
\right)
\right)
\circ^{1\mathbf{Pth}_{\boldsymbol{\mathcal{A}}^{(2)}}}_{s}
\mathrm{ip}^{(2,X)@}_{s}\left(
P
\right)
\right)
\biggr\rrbracket_{s}
\tag{4}
\\&=
\biggl\llbracket 
\mathrm{CH}^{(2)}_{s}\left(
\mathrm{ip}^{(2,X)@}_{s}\left(
\left(
R
\circ^{1\mathbf{PT}_{\boldsymbol{\mathcal{A}}^{(2)}}}_{s}
Q
\right)
\circ^{1\mathbf{PT}_{\boldsymbol{\mathcal{A}}^{(2)}}}_{s}
P
\right)\right)
\biggr\rrbracket_{s}
\tag{5}
\\&=
\biggl\llbracket 
\left(
R
\circ^{1\mathbf{PT}_{\boldsymbol{\mathcal{A}}^{(2)}}}_{s}
Q
\right)
\circ^{1\mathbf{PT}_{\boldsymbol{\mathcal{A}}^{(2)}}}_{s}
P
\biggr\rrbracket_{s}
\tag{6}
\\&=
\left(
\llbracket R\rrbracket_{s}
\circ^{1\llbracket\mathbf{PT}_{\boldsymbol{\mathcal{A}}^{(2)}}\rrbracket}_{s}
\llbracket Q\rrbracket_{s}
\right)
\circ^{1\llbracket\mathbf{PT}_{\boldsymbol{\mathcal{A}}^{(2)}}\rrbracket}_{s}
\llbracket P\rrbracket_{s}.
\tag{7}
\end{align*}
\end{flushleft}

In the just stated chain of equalities, the first equality unravels the description of the operation of $1$-composition in the many sorted partial $\Sigma^{\boldsymbol{\mathcal{A}}^{(2)}}$-algebra $\llbracket\mathbf{PT}_{\boldsymbol{\mathcal{A}}^{(2)}}\rrbracket$ according to Proposition~\ref{PDPTQDCatAlg}; the second equality follows from Proposition~\ref{LDVWCong}; the third equality follows from the fact that $\mathrm{ip}^{(2,X)@}$ is a $\Sigma^{\boldsymbol{\mathcal{A}}^{(2)}}$-homomorphism according to Definition~\ref{DDIp}; the fourth equality follows from Proposition~\ref{PDVDCH} and Proposition~\ref{PDVVarB6}; the fifth equality follows from the fact that $\mathrm{ip}^{(2,X)@}$ is a $\Sigma^{\boldsymbol{\mathcal{A}}^{(2)}}$-homomorphism according to Definition~\ref{DDIp}; the sixth equality follows from Proposition~\ref{LDVWCong}; finally, the last equality recovers the description of the operation of $1$-composition in the many sorted partial $\Sigma^{\boldsymbol{\mathcal{A}}^{(2)}}$-algebra $\llbracket\mathbf{PT}_{\boldsymbol{\mathcal{A}}^{(2)}}\rrbracket$ according to Proposition~\ref{PDPTQDCatAlg}.

This completes the proof.
\end{proof}

\begin{proposition}\label{PDPTQVarAB1} Let $s$ be a sort in $S$ and $\llbracket P\rrbracket_{s}$ a second-order path term class in  $\llbracket\mathrm{PT}_{\boldsymbol{\mathcal{A}}^{(2)}}\rrbracket_{s}$, then the following equalities hold
\allowdisplaybreaks
\begin{align*}
\mathrm{sc}^{1\llbracket\mathbf{PT}_{\boldsymbol{\mathcal{A}}^{(2)}}\rrbracket}_{s}\left(
\mathrm{sc}^{0\llbracket\mathbf{PT}_{\boldsymbol{\mathcal{A}}^{(2)}}\rrbracket}_{s}\left(
\llbracket P\rrbracket_{s}
\right)
\right)
& 
=
\mathrm{sc}^{0\llbracket\mathbf{PT}_{\boldsymbol{\mathcal{A}}^{(2)}}\rrbracket}_{s}\left(
\llbracket P\rrbracket_{s}
\right);
\\
\mathrm{sc}^{0\llbracket\mathbf{PT}_{\boldsymbol{\mathcal{A}}^{(2)}}\rrbracket}_{s}\left(
\mathrm{sc}^{1\llbracket\mathbf{PT}_{\boldsymbol{\mathcal{A}}^{(2)}}\rrbracket}_{s}\left(
\llbracket P\rrbracket_{s}
\right)
\right)
& 
=
\mathrm{sc}^{0\llbracket\mathbf{PT}_{\boldsymbol{\mathcal{A}}^{(2)}}\rrbracket}_{s}\left(
\llbracket P\rrbracket_{s}
\right);
\\
\mathrm{sc}^{0\llbracket\mathbf{PT}_{\boldsymbol{\mathcal{A}}^{(2)}}\rrbracket}_{s}\left(
\mathrm{tg}^{1\llbracket\mathbf{PT}_{\boldsymbol{\mathcal{A}}^{(2)}}\rrbracket}_{s}\left(
\llbracket P\rrbracket_{s}
\right)
\right)
& 
=
\mathrm{sc}^{0\llbracket\mathbf{PT}_{\boldsymbol{\mathcal{A}}^{(2)}}\rrbracket}_{s}\left(
\llbracket P\rrbracket_{s}
\right);
\\
\mathrm{tg}^{1\llbracket\mathbf{PT}_{\boldsymbol{\mathcal{A}}^{(2)}}\rrbracket}_{s}\left(
\mathrm{tg}^{0\llbracket\mathbf{PT}_{\boldsymbol{\mathcal{A}}^{(2)}}\rrbracket}_{s}\left(
\llbracket P\rrbracket_{s}
\right)
\right)
& 
=
\mathrm{tg}^{0\llbracket\mathbf{PT}_{\boldsymbol{\mathcal{A}}^{(2)}}\rrbracket}_{s}\left(
\llbracket P\rrbracket_{s}
\right);
\\
\mathrm{tg}^{0\llbracket\mathbf{PT}_{\boldsymbol{\mathcal{A}}^{(2)}}\rrbracket}_{s}\left(
\mathrm{tg}^{1\llbracket\mathbf{PT}_{\boldsymbol{\mathcal{A}}^{(2)}}\rrbracket}_{s}\left(
\llbracket P\rrbracket_{s}
\right)
\right)
& 
=
\mathrm{tg}^{0\llbracket\mathbf{PT}_{\boldsymbol{\mathcal{A}}^{(2)}}\rrbracket}_{s}\left(
\llbracket P\rrbracket_{s}
\right);
\\
\mathrm{tg}^{0\llbracket\mathbf{PT}_{\boldsymbol{\mathcal{A}}^{(2)}}\rrbracket}_{s}\left(
\mathrm{sc}^{1\llbracket\mathbf{PT}_{\boldsymbol{\mathcal{A}}^{(2)}}\rrbracket}_{s}\left(
\llbracket P\rrbracket_{s}
\right)
\right)
& 
=
\mathrm{tg}^{0\llbracket\mathbf{PT}_{\boldsymbol{\mathcal{A}}^{(2)}}\rrbracket}_{s}\left(
\llbracket P\rrbracket_{s}
\right).
\end{align*}	
\end{proposition}
\begin{proof}
We will only present the proof for the first equality. The other three are handled in a similar manner.

Note that the following chain of equalities holds
\begin{flushleft}
$\mathrm{sc}^{1\llbracket\mathbf{PT}_{\boldsymbol{\mathcal{A}}^{(2)}}\rrbracket}_{s}\left(
\mathrm{sc}^{0\llbracket\mathbf{PT}_{\boldsymbol{\mathcal{A}}^{(2)}}\rrbracket}_{s}\left(
\llbracket P\rrbracket_{s}
\right)
\right)$
\allowdisplaybreaks
\begin{align*}
\quad&=
\biggl\llbracket
\mathrm{sc}^{1\mathbf{PT}_{\boldsymbol{\mathcal{A}}^{(2)}}}_{s}\left(
\mathrm{sc}^{0\mathbf{PT}_{\boldsymbol{\mathcal{A}}^{(2)}}}_{s}\left(
P
\right)\right)
\biggr\rrbracket_{s}
\tag{1}
\\&=
\biggl\llbracket
\mathrm{CH}^{(2)}_{s}\left(
\mathrm{ip}^{(2,X)@}_{s}\left(
\mathrm{sc}^{1\mathbf{PT}_{\boldsymbol{\mathcal{A}}^{(2)}}}_{s}\left(
\mathrm{sc}^{0\mathbf{PT}_{\boldsymbol{\mathcal{A}}^{(2)}}}_{s}\left(
P
\right)\right)
\right)\right)
\biggr\rrbracket_{s}
\tag{2}
\\&=
\biggl\llbracket
\mathrm{CH}^{(2)}_{s}\left(
\mathrm{sc}^{1\mathbf{Pth}_{\boldsymbol{\mathcal{A}}^{(2)}}}_{s}\left(
\mathrm{sc}^{0\mathbf{Pth}_{\boldsymbol{\mathcal{A}}^{(2)}}}_{s}\left(
\mathrm{ip}^{(2,X)@}_{s}\left(
P
\right)\right)
\right)\right)
\biggr\rrbracket_{s}
\tag{3}
\\&=
\biggl\llbracket
\mathrm{CH}^{(2)}_{s}\left(
\mathrm{sc}^{0\mathbf{Pth}_{\boldsymbol{\mathcal{A}}^{(2)}}}_{s}\left(
\mathrm{ip}^{(2,X)@}_{s}\left(
P
\right)
\right)\right)
\biggr\rrbracket_{s}
\tag{4}
\\&=
\biggl\llbracket
\mathrm{CH}^{(2)}_{s}\left(
\mathrm{ip}^{(2,X)@}_{s}\left(
\mathrm{sc}^{0\mathbf{PT}_{\boldsymbol{\mathcal{A}}^{(2)}}}_{s}\left(
P
\right)
\right)\right)
\biggr\rrbracket_{s}
\tag{5}
\\&=
\biggl\llbracket
\mathrm{sc}^{0\mathbf{PT}_{\boldsymbol{\mathcal{A}}^{(2)}}}_{s}\left(
P
\right)
\biggr\rrbracket_{s}
\tag{6}
\\&=
\mathrm{sc}^{0\llbracket\mathbf{PT}_{\boldsymbol{\mathcal{A}}^{(2)}}\rrbracket}_{s}\left(
\llbracket P\rrbracket_{s}
\right).
\tag{7}
\end{align*}
\end{flushleft}

In the just stated chain of equalities, the first equality unravels the description of the operations of $1$-source and $0$-source in the many sorted partial $\Sigma^{\boldsymbol{\mathcal{A}}^{(2)}}$-algebra $\llbracket\mathbf{PT}_{\boldsymbol{\mathcal{A}}^{(2)}}\rrbracket$ according to Proposition~\ref{PDPTQDCatAlg}; the second equality follows from Proposition~\ref{LDVWCong}; the third equality follows from the fact that $\mathrm{ip}^{(2,X)@}$ is a $\Sigma^{\boldsymbol{\mathcal{A}}^{(2)}}$-homomorphism according to Definition~\ref{DDIp}; the fourth equality follows from Proposition~\ref{PDVDCH} and Proposition~\ref{PDVVarAB1}; the fifth equality follows from the fact that $\mathrm{ip}^{(2,X)@}$ is a $\Sigma^{\boldsymbol{\mathcal{A}}^{(2)}}$-homomorphism according to Definition~\ref{DDIp}; the sixth equality follows from Proposition~\ref{LDVWCong}; finally, the last equality recovers the description of the operation of $0$-source in the many sorted partial $\Sigma^{\boldsymbol{\mathcal{A}}^{(2)}}$-algebra $\llbracket\mathbf{PT}_{\boldsymbol{\mathcal{A}}^{(2)}}\rrbracket$ according to Proposition~\ref{PDPTQDCatAlg}.

This completes the proof.
\end{proof}

\begin{proposition}\label{PDPTQVarAB2} Let $s$ be a sort in $S$ and $\llbracket P\rrbracket_{s}$, $\llbracket Q\rrbracket_{s}$ second-order path term classes in  $\llbracket\mathrm{PT}_{\boldsymbol{\mathcal{A}}^{(2)}}\rrbracket_{s}$ satisfying that
\[\mathrm{sc}^{0\llbracket\mathbf{PT}_{\boldsymbol{\mathcal{A}}^{(2)}}\rrbracket}_{s}\left(
\llbracket Q\rrbracket_{s}
\right)=
\mathrm{tg}^{0\llbracket\mathbf{PT}_{\boldsymbol{\mathcal{A}}^{(2)}}\rrbracket}_{s}\left(
\llbracket P\rrbracket_{s}
\right),
\]
then the following equalities hold
\allowdisplaybreaks
\begin{multline*}
\mathrm{sc}^{1\llbracket\mathbf{PT}_{\boldsymbol{\mathcal{A}}^{(2)}}\rrbracket}_{s}\left(
\llbracket Q\rrbracket_{s}
\circ^{0\llbracket\mathbf{PT}_{\boldsymbol{\mathcal{A}}^{(2)}}\rrbracket}_{s}
\llbracket P\rrbracket_{s}
\right)
\\=
\mathrm{sc}^{1\llbracket\mathbf{PT}_{\boldsymbol{\mathcal{A}}^{(2)}}\rrbracket}_{s}\left(
\llbracket Q\rrbracket_{s}
\right)
\circ^{0\llbracket\mathbf{PT}_{\boldsymbol{\mathcal{A}}^{(2)}}\rrbracket}_{s}
\mathrm{sc}^{1\llbracket\mathbf{PT}_{\boldsymbol{\mathcal{A}}^{(2)}}\rrbracket}_{s}\left(
\llbracket P\rrbracket_{s}
\right);
\end{multline*}
\allowdisplaybreaks
\begin{multline*}
\mathrm{tg}^{1\llbracket\mathbf{PT}_{\boldsymbol{\mathcal{A}}^{(2)}}\rrbracket}_{s}\left(
\llbracket Q\rrbracket_{s}
\circ^{0\llbracket\mathbf{PT}_{\boldsymbol{\mathcal{A}}^{(2)}}\rrbracket}_{s}
\llbracket P\rrbracket_{s}
\right)
\\=
\mathrm{tg}^{1\llbracket\mathbf{PT}_{\boldsymbol{\mathcal{A}}^{(2)}}\rrbracket}_{s}\left(
\llbracket Q\rrbracket_{s}
\right)
\circ^{0\llbracket\mathbf{PT}_{\boldsymbol{\mathcal{A}}^{(2)}}\rrbracket}_{s}
\mathrm{tg}^{1\llbracket\mathbf{PT}_{\boldsymbol{\mathcal{A}}^{(2)}}\rrbracket}_{s}\left(
\llbracket P\rrbracket_{s}
\right).
\end{multline*}
\end{proposition}
\begin{proof}

We will only present the proof of the first equality. The other one is handled in a similar manner.

The following chain of equalities holds
\begin{flushleft}
$\mathrm{sc}^{1\llbracket\mathbf{PT}_{\boldsymbol{\mathcal{A}}^{(2)}}\rrbracket}_{s}\left(
\llbracket Q\rrbracket_{s}
\circ^{0\llbracket\mathbf{PT}_{\boldsymbol{\mathcal{A}}^{(2)}}\rrbracket}_{s}
\llbracket P\rrbracket_{s}
\right)$
\allowdisplaybreaks
\begin{align*}
&=
\biggl\llbracket 
\mathrm{sc}^{1\mathbf{PT}_{\boldsymbol{\mathcal{A}}^{(2)}}}_{s}\left(
Q
\circ^{0\mathbf{PT}_{\boldsymbol{\mathcal{A}}^{(2)}}}_{s}
P
\right)
\biggr\rrbracket_{s}
\tag{1}
\\&=
\biggl\llbracket 
\mathrm{CH}^{(2)}_{s}\left(
\mathrm{ip}^{(2,X)@}_{s}\left(
\mathrm{sc}^{1\mathbf{PT}_{\boldsymbol{\mathcal{A}}^{(2)}}}_{s}\left(
Q
\circ^{0\mathbf{PT}_{\boldsymbol{\mathcal{A}}^{(2)}}}_{s}
P
\right)\right)\right)
\biggr\rrbracket_{s}
\tag{2}
\\&=
\biggl\llbracket 
\mathrm{CH}^{(2)}_{s}\left(
\mathrm{sc}^{1\mathbf{Pth}_{\boldsymbol{\mathcal{A}}^{(2)}}}_{s}\left(
\mathrm{ip}^{(2,X)@}_{s}\left(
Q
\right)
\circ^{0\mathbf{Pth}_{\boldsymbol{\mathcal{A}}^{(2)}}}_{s}
\mathrm{ip}^{(2,X)@}_{s}\left(
P
\right)
\right)
\right)
\biggr\rrbracket_{s}
\tag{3}
\\&=
\biggl\llbracket 
\mathrm{CH}^{(2)}_{s}\left(
\mathrm{sc}^{1\mathbf{Pth}_{\boldsymbol{\mathcal{A}}^{(2)}}}_{s}\left(
\mathrm{ip}^{(2,X)@}_{s}\left(
Q
\right)
\right)
\circ^{0\mathbf{Pth}_{\boldsymbol{\mathcal{A}}^{(2)}}}_{s}
\mathrm{sc}^{1\mathbf{Pth}_{\boldsymbol{\mathcal{A}}^{(2)}}}_{s}\left(
\mathrm{ip}^{(2,X)@}_{s}\left(
P
\right)
\right)
\right)
\biggr\rrbracket_{s}
\tag{4}
\\&=
\biggl\llbracket 
\mathrm{CH}^{(2)}_{s}\left(
\mathrm{ip}^{(2,X)@}_{s}\left(
\mathrm{sc}^{1\mathbf{PT}_{\boldsymbol{\mathcal{A}}^{(2)}}}_{s}\left(
Q
\right)
\circ^{0\mathbf{PT}_{\boldsymbol{\mathcal{A}}^{(2)}}}_{s}
\mathrm{sc}^{1\mathbf{PT}_{\boldsymbol{\mathcal{A}}^{(2)}}}_{s}\left(
P
\right)
\right)\right)
\biggr\rrbracket_{s}
\tag{5}
\\&=
\biggl\llbracket 
\mathrm{sc}^{1\mathbf{PT}_{\boldsymbol{\mathcal{A}}^{(2)}}}_{s}\left(
Q
\right)
\circ^{0\mathbf{PT}_{\boldsymbol{\mathcal{A}}^{(2)}}}_{s}
\mathrm{sc}^{1\mathbf{PT}_{\boldsymbol{\mathcal{A}}^{(2)}}}_{s}\left(
P
\right)
\biggr\rrbracket_{s}
\tag{6}
\\&=
\mathrm{sc}^{1\llbracket\mathbf{PT}_{\boldsymbol{\mathcal{A}}^{(2)}}\rrbracket}_{s}\left(
\llbracket Q\rrbracket_{s}
\right)
\circ^{0\llbracket\mathbf{PT}_{\boldsymbol{\mathcal{A}}^{(2)}}\rrbracket}_{s}
\mathrm{sc}^{1\llbracket\mathbf{PT}_{\boldsymbol{\mathcal{A}}^{(2)}}\rrbracket}_{s}\left(
\llbracket P\rrbracket_{s}
\right).
\tag{7}
\end{align*}
\end{flushleft}

In the just stated chain of equalities, the first equality unravels the description of the operation of $1$-source and $0$-composition in the many sorted partial $\Sigma^{\boldsymbol{\mathcal{A}}^{(2)}}$-algebra $\llbracket\mathbf{PT}_{\boldsymbol{\mathcal{A}}^{(2)}}\rrbracket$ according to Proposition~\ref{PDPTQDCatAlg}; the second equality follows from Proposition~\ref{LDVWCong}; the third equality follows from the fact that $\mathrm{ip}^{(2,X)@}$ is a $\Sigma^{\boldsymbol{\mathcal{A}}^{(2)}}$-homomorphism according to Definition~\ref{DDIp}; the fourth equality follows from Proposition~\ref{PDVDCH} and Proposition~\ref{PDVVarAB2}; the fifth equality follows from the fact that $\mathrm{ip}^{(2,X)@}$ is a $\Sigma^{\boldsymbol{\mathcal{A}}^{(2)}}$-homomorphism according to Definition~\ref{DDIp}; the sixth equality follows from Proposition~\ref{LDVWCong}; finally, the last equality recovers the description of the operations of $1$-source nd $0$-composition in the many sorted partial $\Sigma^{\boldsymbol{\mathcal{A}}^{(2)}}$-algebra $\llbracket\mathbf{PT}_{\boldsymbol{\mathcal{A}}^{(2)}}\rrbracket$ according to Proposition~\ref{PDPTQDCatAlg}.

This completes the proof.
\end{proof}

\begin{proposition}\label{PDPTQVarAB3} Let $s$ be a sort in $S$ and $\llbracket P\rrbracket_{s}$, $\llbracket P'\rrbracket_{s}$, $\llbracket Q\rrbracket_{s}$, $\llbracket Q'\rrbracket_{s}$ be second-order path term classes in  $\llbracket\mathrm{PT}_{\boldsymbol{\mathcal{A}}^{(2)}}\rrbracket_{s}$ satisfying that
\allowdisplaybreaks
\begin{align*}
\mathrm{sc}^{1\llbracket\mathbf{PT}_{\boldsymbol{\mathcal{A}}^{(2)}}\rrbracket}_{s}\left(
\llbracket Q'\rrbracket_{s}
\right)
&=
\mathrm{tg}^{1\llbracket\mathbf{PT}_{\boldsymbol{\mathcal{A}}^{(2)}}\rrbracket}_{s}\left(
\llbracket Q\rrbracket_{s}
\right);
\\
\mathrm{sc}^{1\llbracket\mathbf{PT}_{\boldsymbol{\mathcal{A}}^{(2)}}\rrbracket}_{s}\left(
\llbracket P'\rrbracket_{s}
\right)
&=
\mathrm{tg}^{1\llbracket\mathbf{PT}_{\boldsymbol{\mathcal{A}}^{(2)}}\rrbracket}_{s}\left(
\llbracket P\rrbracket_{s}
\right);
\\
\mathrm{sc}^{0\llbracket\mathbf{PT}_{\boldsymbol{\mathcal{A}}^{(2)}}\rrbracket}_{s}\left(
\llbracket Q'\rrbracket_{s}
\right)
&=
\mathrm{tg}^{0\llbracket\mathbf{PT}_{\boldsymbol{\mathcal{A}}^{(2)}}\rrbracket}_{s}\left(
\llbracket P'\rrbracket_{s}
\right);
\\
\mathrm{sc}^{0\llbracket\mathbf{PT}_{\boldsymbol{\mathcal{A}}^{(2)}}\rrbracket}_{s}\left(
\llbracket Q\rrbracket_{s}
\right)
&=
\mathrm{tg}^{0\llbracket\mathbf{PT}_{\boldsymbol{\mathcal{A}}^{(2)}}\rrbracket}_{s}\left(
\llbracket P\rrbracket_{s}
\right),
\end{align*}
then the following equality holds
\allowdisplaybreaks
\begin{multline*}
\left(
\llbracket Q'\rrbracket_{s}
\circ^{0\llbracket\mathbf{PT}_{\boldsymbol{\mathcal{A}}^{(2)}}\rrbracket}_{s}
\llbracket P'\rrbracket_{s}
\right)
\circ^{1\llbracket\mathbf{PT}_{\boldsymbol{\mathcal{A}}^{(2)}}\rrbracket}_{s}
\left(
\llbracket Q\rrbracket_{s}
\circ^{0\llbracket\mathbf{PT}_{\boldsymbol{\mathcal{A}}^{(2)}}\rrbracket}_{s}
\llbracket P\rrbracket_{s}
\right)
\\=
\left(
\llbracket Q'\rrbracket_{s}
\circ^{1\llbracket\mathbf{PT}_{\boldsymbol{\mathcal{A}}^{(2)}}\rrbracket}_{s}
\llbracket Q\rrbracket_{s}
\right)
\circ^{0\llbracket\mathbf{PT}_{\boldsymbol{\mathcal{A}}^{(2)}}\rrbracket}_{s}
\left(
\llbracket P'\rrbracket_{s}
\circ^{1\llbracket\mathbf{PT}_{\boldsymbol{\mathcal{A}}^{(2)}}\rrbracket}_{s}
\llbracket P\rrbracket_{s}
\right).
\end{multline*}
\end{proposition}
\begin{proof}

The following chain of equalities holds
\begin{flushleft}
$
\left(
\llbracket Q'\rrbracket_{s}
\circ^{0\llbracket\mathbf{PT}_{\boldsymbol{\mathcal{A}}^{(2)}}\rrbracket}_{s}
\llbracket P'\rrbracket_{s}
\right)
\circ^{1\llbracket\mathbf{PT}_{\boldsymbol{\mathcal{A}}^{(2)}}\rrbracket}_{s}
\left(
\llbracket Q\rrbracket_{s}
\circ^{0\llbracket\mathbf{PT}_{\boldsymbol{\mathcal{A}}^{(2)}}\rrbracket}_{s}
\llbracket P\rrbracket_{s}
\right)$
\allowdisplaybreaks
\begin{align*}
&=
\biggl\llbracket
\left(
Q'
\circ^{0\mathbf{PT}_{\boldsymbol{\mathcal{A}}^{(2)}}}_{s}
P'
\right)
\circ^{1\mathbf{PT}_{\boldsymbol{\mathcal{A}}^{(2)}}}_{s}
\left(
Q
\circ^{0\mathbf{PT}_{\boldsymbol{\mathcal{A}}^{(2)}}}_{s}
P
\right)
\biggr\rrbracket_{s}
\tag{1}
\\&=
\biggl\llbracket
\mathrm{CH}^{(2)}_{s}\left(
\mathrm{ip}^{(2,X)@}_{s}\left(
\left(
Q'
\circ^{0\mathbf{PT}_{\boldsymbol{\mathcal{A}}^{(2)}}}_{s}
P'
\right)
\circ^{1\mathbf{PT}_{\boldsymbol{\mathcal{A}}^{(2)}}}_{s}
\left(
Q
\circ^{0\mathbf{PT}_{\boldsymbol{\mathcal{A}}^{(2)}}}_{s}
P
\right)
\right)\right)
\biggr\rrbracket_{s}
\tag{2}
\\&=
\biggl\llbracket
\mathrm{CH}^{(2)}_{s}\left(
\left(
\mathrm{ip}^{(2,X)@}_{s}\left(
Q'
\right)
\circ^{0\mathbf{Pth}_{\boldsymbol{\mathcal{A}}^{(2)}}}_{s}
\mathrm{ip}^{(2,X)@}_{s}\left(
P'
\right)
\right)
\circ^{1\mathbf{Pth}_{\boldsymbol{\mathcal{A}}^{(2)}}}_{s}
\right.
\\&\qquad\qquad\qquad\qquad\qquad\qquad\qquad
\left.
\left(
\mathrm{ip}^{(2,X)@}_{s}\left(
Q
\right)
\circ^{0\mathbf{Pth}_{\boldsymbol{\mathcal{A}}^{(2)}}}_{s}
\mathrm{ip}^{(2,X)@}_{s}\left(
P
\right)
\right)\right)
\biggr\rrbracket_{s}
\tag{3}
\\&=
\biggl\llbracket
\mathrm{CH}^{(2)}_{s}\left(
\left(
\mathrm{ip}^{(2,X)@}_{s}\left(
Q'
\right)
\circ^{1\mathbf{Pth}_{\boldsymbol{\mathcal{A}}^{(2)}}}_{s}
\mathrm{ip}^{(2,X)@}_{s}\left(
Q
\right)
\right)
\circ^{0\mathbf{Pth}_{\boldsymbol{\mathcal{A}}^{(2)}}}_{s}
\right.
\\&\qquad\qquad\qquad\qquad\qquad\qquad\qquad
\left.
\left(
\mathrm{ip}^{(2,X)@}_{s}\left(
P'
\right)
\circ^{1\mathbf{Pth}_{\boldsymbol{\mathcal{A}}^{(2)}}}_{s}
\mathrm{ip}^{(2,X)@}_{s}\left(
P
\right)
\right)\right)
\biggr\rrbracket_{s}
\tag{4}
\\&=
\biggl\llbracket
\mathrm{CH}^{(2)}_{s}\left(
\mathrm{ip}^{(2,X)@}_{s}\left(
\left(
Q'
\circ^{1\mathbf{PT}_{\boldsymbol{\mathcal{A}}^{(2)}}}_{s}
Q
\right)
\circ^{0\mathbf{PT}_{\boldsymbol{\mathcal{A}}^{(2)}}}_{s}
\left(
P'
\circ^{1\mathbf{PT}_{\boldsymbol{\mathcal{A}}^{(2)}}}_{s}
P
\right)
\right)\right)
\biggr\rrbracket_{s}
\tag{5}
\\&=
\biggl\llbracket
\left(
Q'
\circ^{1\mathbf{PT}_{\boldsymbol{\mathcal{A}}^{(2)}}}_{s}
Q
\right)
\circ^{0\mathbf{PT}_{\boldsymbol{\mathcal{A}}^{(2)}}}_{s}
\left(
P'
\circ^{1\mathbf{PT}_{\boldsymbol{\mathcal{A}}^{(2)}}}_{s}
P
\right)
\biggr\rrbracket_{s}
\tag{6}
\\&=
\left(
\llbracket Q'\rrbracket_{s}
\circ^{1\llbracket\mathbf{PT}_{\boldsymbol{\mathcal{A}}^{(2)}}\rrbracket}_{s}
\llbracket Q\rrbracket_{s}
\right)
\circ^{0\llbracket\mathbf{PT}_{\boldsymbol{\mathcal{A}}^{(2)}}\rrbracket}_{s}
\left(
\llbracket P'\rrbracket_{s}
\circ^{1\llbracket\mathbf{PT}_{\boldsymbol{\mathcal{A}}^{(2)}}\rrbracket}_{s}
\llbracket P\rrbracket_{s}
\right).
\tag{7}
\end{align*}
\end{flushleft}

In the just stated chain of equalities, the first equality unravels the description of the operations of $0$-composition and $1$-composition in the many sorted partial $\Sigma^{\boldsymbol{\mathcal{A}}^{(2)}}$-algebra $\llbracket\mathbf{PT}_{\boldsymbol{\mathcal{A}}^{(2)}}\rrbracket$ according to Proposition~\ref{PDPTQDCatAlg}; the second equality follows from Proposition~\ref{LDVWCong}; the third equality follows from the fact that $\mathrm{ip}^{(2,X)@}$ is a $\Sigma^{\boldsymbol{\mathcal{A}}^{(2)}}$-homomorphism according to Definition~\ref{DDIp}; the fourth equality follows from Proposition~\ref{PDVDCH} and Proposition~\ref{PDVVarAB3}; the fifth equality follows from the fact that $\mathrm{ip}^{(2,X)@}$ is a $\Sigma^{\boldsymbol{\mathcal{A}}^{(2)}}$-homomorphism according to Definition~\ref{DDIp}; the sixth equality follows from Proposition~\ref{LDVWCong}; finally, the last equality recovers the description of the operations of $0$-composition and $1$-composition in the many sorted partial $\Sigma^{\boldsymbol{\mathcal{A}}^{(2)}}$-algebra $\llbracket\mathbf{PT}_{\boldsymbol{\mathcal{A}}^{(2)}}\rrbracket$ according to Proposition~\ref{PDPTQDCatAlg}.

This completes the proof.
\end{proof}

In virtue of the foregoing propositions, we can consider the $S$-sorted category given by second-order path term $\Theta^{\llbracket 2\rrbracket}$-classes.

\begin{restatable}{definition}{DDPTQDCat}
\label{DDPTQDCat}
\index{path terms!second-order!$\llbracket\mathsf{PT}_{\boldsymbol{\mathcal{A}}^{(2)}}\rrbracket$}
Let $\llbracket\mathsf{PT}_{\boldsymbol{\mathcal{A}}^{(2)}}\rrbracket$ denote the ordered tuple
\allowdisplaybreaks
\begin{multline*}
\llbracket\mathsf{PT}_{\boldsymbol{\mathcal{A}}^{(2)}}\rrbracket
=
\biggl(
\llbracket\mathbf{PT}_{\boldsymbol{\mathcal{A}}^{(2)}}\rrbracket,
\\
\left(
\circ^{0\llbracket\mathbf{PT}_{\boldsymbol{\mathcal{A}}^{(2)}}\rrbracket},
\mathrm{sc}^{0\llbracket\mathbf{PT}_{\boldsymbol{\mathcal{A}}^{(2)}}\rrbracket},
\mathrm{tg}^{0\llbracket\mathbf{PT}_{\boldsymbol{\mathcal{A}}^{(2)}}\rrbracket}
\right),
\\
\left.
\left(
\circ^{1\llbracket\mathbf{PT}_{\boldsymbol{\mathcal{A}}^{(2)}}\rrbracket},
\mathrm{sc}^{1\llbracket\mathbf{PT}_{\boldsymbol{\mathcal{A}}^{(2)}}\rrbracket},
\mathrm{tg}^{1\llbracket\mathbf{PT}_{\boldsymbol{\mathcal{A}}^{(2)}}\rrbracket}
\right)
\right).
\end{multline*}
\end{restatable}

\begin{restatable}{proposition}{PDPTQDCat}
\label{PDPTQDCat} $\llbracket\mathsf{PT}_{\boldsymbol{\mathcal{A}}^{(2)}}\rrbracket$ is an $S$-sorted $2$-category.
\end{restatable}
\begin{proof}
Let us recall that, for every sort $s\in S$, the structure 
$(
\llbracket\mathbf{PT}_{\boldsymbol{\mathcal{A}}^{(2)}}\rrbracket_{s}, 
\xi_{0,s},
\xi_{1,s}
)$, where
\begin{align*}
\xi_{0,s}&=
\left(
\circ^{0\llbracket\mathbf{PT}_{\boldsymbol{\mathcal{A}}^{(2)}}\rrbracket}_{s},
\mathrm{sc}^{0\llbracket\mathbf{PT}_{\boldsymbol{\mathcal{A}}^{(2)}}\rrbracket}_{s},
\mathrm{tg}^{0\llbracket\mathbf{PT}_{\boldsymbol{\mathcal{A}}^{(2)}}\rrbracket}_{s}
\right);
\\
\xi_{1,s}&=
\left(
\circ^{1\llbracket\mathbf{PT}_{\boldsymbol{\mathcal{A}}^{(2)}}\rrbracket}_{s},
\mathrm{sc}^{1\llbracket\mathbf{PT}_{\boldsymbol{\mathcal{A}}^{(2)}}\rrbracket}_{s},
\mathrm{tg}^{1\llbracket\mathbf{PT}_{\boldsymbol{\mathcal{A}}^{(2)}}\rrbracket}_{s}
\right),
\end{align*}
is a single-sorted  $2$-category in virtue of Definition~\ref{D2Cat} and Propositions~\ref{PDPTQVarA2}, \ref{PDPTQVarA3}, \ref{PDPTQVarA4}, \ref{PDPTQVarA5}, \ref{PDPTQVarA6}, \ref{PDPTQVarB2}, \ref{PDPTQVarB3}, \ref{PDPTQVarB4}, \ref{PDPTQVarB5}, \ref{PDPTQVarB6}, \ref{PDPTQVarAB1}, \ref{PDPTQVarAB2}, and \ref{PDPTQVarAB3}. Therefore, following Definition~\ref{DnCat}, the $S$-sorted structure 
\[
\left(
\llbracket\mathbf{PT}_{\boldsymbol{\mathcal{A}}^{(2)}}\rrbracket_{s},
\xi_{0,s},
\xi_{1,s}
\right)_{s\in S}
\]   
is an $S$-sorted $2$-category.
\end{proof}

The most important feature, as we will see immediately, is that $\llbracket\mathsf{PT}_{\boldsymbol{\mathcal{A}}^{(2)}}\rrbracket$ is an $S$-sorted $2$-categorial $\Sigma$-algebra. To state this fact properly, we present the following propositions.

\begin{proposition}\label{PDPTQVarA7} Let $(\mathbf{s},s)$ be an element in $S^{\star}\times S$, $\sigma$ an operation symbol in $\Sigma_{\mathbf{s},s}$ and 
$(\llbracket P_{j}\rrbracket_{s_{j}})_{j\in\bb{\mathbf{s}}}$ a family of second-order path term classes in  $\llbracket\mathrm{PT}_{\boldsymbol{\mathcal{A}}^{(2)}}\rrbracket_{\mathbf{s}}$, then the following equalities hold
\allowdisplaybreaks
\begin{align*}
\sigma^{\llbracket\mathbf{PT}_{\boldsymbol{\mathcal{A}}^{(2)}}\rrbracket}
\left(\left(
\mathrm{sc}^{0\llbracket\mathbf{PT}_{\boldsymbol{\mathcal{A}}^{(2)}}\rrbracket}_{s_{j}}
\left(
\llbracket P_{j}\rrbracket_{s_{j}}
\right)\right)_{j\in\bb{\mathbf{s}}}
\right)
&=
\mathrm{sc}^{0\llbracket\mathbf{PT}_{\boldsymbol{\mathcal{A}}^{(2)}}\rrbracket}_{s}\left(
\sigma^{\llbracket\mathbf{PT}_{\boldsymbol{\mathcal{A}}^{(2)}}\rrbracket}
\left(\left(
\llbracket P_{j}\rrbracket_{s_{j}}
\right)_{j\in\bb{\mathbf{s}}}
\right)
\right);
\\
\sigma^{\llbracket\mathbf{PT}_{\boldsymbol{\mathcal{A}}^{(2)}}\rrbracket}
\left(\left(
\mathrm{tg}^{0\llbracket\mathbf{PT}_{\boldsymbol{\mathcal{A}}^{(2)}}\rrbracket}_{s_{j}}
\left(
\llbracket P_{j}\rrbracket_{s_{j}}
\right)\right)_{j\in\bb{\mathbf{s}}}
\right)
&=
\mathrm{tg}^{0\llbracket\mathbf{PT}_{\boldsymbol{\mathcal{A}}^{(2)}}\rrbracket}_{s}\left(
\sigma^{\llbracket\mathbf{PT}_{\boldsymbol{\mathcal{A}}^{(2)}}\rrbracket}
\left(\left(
\llbracket P_{j}\rrbracket_{s_{j}}
\right)_{j\in\bb{\mathbf{s}}}
\right)
\right).
\end{align*}
\end{proposition}
\begin{proof}
We will only present the proof for the first equality. The other one is handled in a similar manner.

Note that the following chain of equalities holds
\begin{flushleft}
$\sigma^{\llbracket\mathbf{PT}_{\boldsymbol{\mathcal{A}}^{(2)}}\rrbracket}
\left(\left(
\mathrm{sc}^{0\llbracket\mathbf{PT}_{\boldsymbol{\mathcal{A}}^{(2)}}\rrbracket}_{s_{j}}
\left(
\llbracket P_{j}\rrbracket_{s_{j}}
\right)\right)_{j\in\bb{\mathbf{s}}}
\right)$
\allowdisplaybreaks
\begin{align*}
&=
\Biggl\llbracket
\sigma^{\mathbf{PT}_{\boldsymbol{\mathcal{A}}^{(2)}}}
\left(\left(
\mathrm{sc}^{0\mathbf{PT}_{\boldsymbol{\mathcal{A}}^{(2)}}}_{s_{j}}
\left(
P_{j}
\right)
\right)_{j\in\bb{\mathbf{s}}}
\right)
\Biggr\rrbracket_{s}
\tag{1}
\\
&=
\Biggl\llbracket
\mathrm{CH}^{(2)}_{s}\left(
\mathrm{ip}^{(2,X)@}_{s}\left(
\sigma^{\mathbf{PT}_{\boldsymbol{\mathcal{A}}^{(2)}}}
\left(\left(
\mathrm{sc}^{0\mathbf{PT}_{\boldsymbol{\mathcal{A}}^{(2)}}}_{s_{j}}
\left(
P_{j}
\right)
\right)_{j\in\bb{\mathbf{s}}}
\right)
\right)
\right)
\Biggr\rrbracket_{s}
\tag{2}
\\
&=
\Biggl\llbracket
\mathrm{CH}^{(2)}_{s}\left(
\sigma^{\mathbf{Pth}_{\boldsymbol{\mathcal{A}}^{(2)}}}
\left(\left(
\mathrm{sc}^{0\mathbf{Pth}_{\boldsymbol{\mathcal{A}}^{(2)}}}_{s_{j}}
\left(
\mathrm{ip}^{(2,X)@}_{s_{j}}\left(
P_{j}
\right)
\right)
\right)_{j\in\bb{\mathbf{s}}}
\right)
\right)
\Biggr\rrbracket_{s}
\tag{3}
\\
&=
\Biggl\llbracket
\mathrm{CH}^{(2)}_{s}\left(
\mathrm{sc}^{0\mathbf{Pth}_{\boldsymbol{\mathcal{A}}^{(2)}}}_{s}
\left(
\sigma^{\mathbf{Pth}_{\boldsymbol{\mathcal{A}}^{(2)}}}
\left(\left(
\mathrm{ip}^{(2,X)@}_{s_{j}}\left(
P_{j}
\right)
\right)_{j\in\bb{\mathbf{s}}}
\right)
\right)
\right)
\Biggr\rrbracket_{s}
\tag{4}
\\
&=
\Biggl\llbracket
\mathrm{CH}^{(2)}_{s}\left(
\mathrm{ip}^{(2,X)@}_{s}\left(
\mathrm{sc}^{0\mathbf{PT}_{\boldsymbol{\mathcal{A}}^{(2)}}}_{s}
\left(
\sigma^{\mathbf{PT}_{\boldsymbol{\mathcal{A}}^{(2)}}}
\left(\left(
P_{j}
\right)_{j\in\bb{\mathbf{s}}}
\right)
\right)
\right)
\right)
\Biggr\rrbracket_{s}
\tag{5}
\\
&=
\Biggl\llbracket
\mathrm{sc}^{0\mathbf{PT}_{\boldsymbol{\mathcal{A}}^{(2)}}}_{s}
\left(
\sigma^{\mathbf{PT}_{\boldsymbol{\mathcal{A}}^{(2)}}}
\left(\left(
P_{j}
\right)_{j\in\bb{\mathbf{s}}}
\right)
\right)
\Biggr\rrbracket_{s}
\tag{6}
\\&=
\mathrm{sc}^{0\llbracket\mathbf{PT}_{\boldsymbol{\mathcal{A}}^{(2)}}\rrbracket}_{s}\left(
\sigma^{\llbracket\mathbf{PT}_{\boldsymbol{\mathcal{A}}^{(2)}}\rrbracket}
\left(\left(
\llbracket P_{j}\rrbracket_{s_{j}}
\right)_{j\in\bb{\mathbf{s}}}
\right)
\right).
\tag{7}
\end{align*}
\end{flushleft}

In the just stated chain of equalities, the first equality unravels the description of the operations of $0$-source and $\sigma$ in the many sorted partial $\Sigma^{\boldsymbol{\mathcal{A}}^{(2)}}$-algebra $\llbracket\mathbf{PT}_{\boldsymbol{\mathcal{A}}^{(2)}}\rrbracket$ according to Proposition~\ref{PDPTQDCatAlg}; the second equality follows from Proposition~\ref{LDVWCong}; the third equality follows from the fact that $\mathrm{ip}^{(2,X)@}$ is a $\Sigma^{\boldsymbol{\mathcal{A}}^{(2)}}$-homomorphism according to Definition~\ref{DDIp}; the fourth equality follows from Proposition~\ref{PDVDCH} and Proposition~\ref{PDVVarA7}; the fifth equality follows from the fact that $\mathrm{ip}^{(2,X)@}$ is a $\Sigma^{\boldsymbol{\mathcal{A}}^{(2)}}$-homomorphism according to Definition~\ref{DDIp}; the sixth equality follows from Proposition~\ref{LDVWCong}; finally, the last equality recovers the description of the operations of $0$-source and $\sigma$ in the many sorted partial $\Sigma^{\boldsymbol{\mathcal{A}}^{(2)}}$-algebra $\llbracket\mathbf{PT}_{\boldsymbol{\mathcal{A}}^{(2)}}\rrbracket$ according to Proposition~\ref{PDPTQDCatAlg}.

This completes the proof.
\end{proof}

\begin{proposition}\label{PDPTQVarA8} Let $(\mathbf{s},s)$ be an element in $S^{\star}\times S$, $\sigma$ and operation symbol in $\Sigma_{\mathbf{s},s}$ and $(\llbracket Q_{j}\rrbracket_{s_{j}})_{j\in\bb{\mathbf{s}}}$, $(\llbracket P_{j}\rrbracket_{s_{j}})_{j\in\bb{\mathbf{s}}}$ be two families of second-order path term classes in  $\llbracket\mathrm{PT}_{\boldsymbol{\mathcal{A}}^{(2)}}\rrbracket_{\mathbf{s}}$ satisfying that, for every $j\in\bb{\mathbf{s}}$,
\[\mathrm{sc}^{0\llbracket\mathbf{PT}_{\boldsymbol{\mathcal{A}}^{(2)}}\rrbracket}_{s_{j}}\left(
\llbracket Q_{j}\rrbracket_{s_{j}}
\right)=
\mathrm{tg}^{0\llbracket\mathbf{PT}_{\boldsymbol{\mathcal{A}}^{(2)}}\rrbracket}_{s_{j}}\left(
\llbracket P_{j}\rrbracket_{s_{j}}
\right),
\]
then the following equality holds
\allowdisplaybreaks
\begin{multline*}
\sigma^{\llbracket\mathbf{PT}_{\boldsymbol{\mathcal{A}}^{(2)}}\rrbracket}
\left(\left(
\llbracket Q_{j}\rrbracket_{s_{j}}
\circ^{0\llbracket\mathbf{PT}_{\boldsymbol{\mathcal{A}}^{(2)}}\rrbracket}_{s_{j}}
\llbracket P_{j}\rrbracket_{s_{j}}
\right)_{j\in\bb{\mathbf{s}}}
\right)
\\=
\sigma^{\llbracket\mathbf{PT}_{\boldsymbol{\mathcal{A}}^{(2)}}\rrbracket}
\left(\left(
\llbracket Q_{j}\rrbracket_{s_{j}}
\right)_{j\in\bb{\mathbf{s}}}\right)
\circ^{0\llbracket\mathbf{PT}_{\boldsymbol{\mathcal{A}}^{(2)}}\rrbracket}_{s}
\sigma^{\llbracket\mathbf{PT}_{\boldsymbol{\mathcal{A}}^{(2)}}\rrbracket}
\left(\left(
\llbracket P_{j}\rrbracket_{s_{j}}
\right)_{j\in\bb{\mathbf{s}}}\right).
\end{multline*}
\end{proposition}
\begin{proof}
Note that the following chain of equalities holds
\begin{flushleft}
$\sigma^{\llbracket\mathbf{PT}_{\boldsymbol{\mathcal{A}}^{(2)}}\rrbracket}
\left(\left(
\llbracket Q_{j}\rrbracket_{s_{j}}
\circ^{0\llbracket\mathbf{PT}_{\boldsymbol{\mathcal{A}}^{(2)}}\rrbracket}_{s_{j}}
\llbracket P_{j}\rrbracket_{s_{j}}
\right)_{j\in\bb{\mathbf{s}}}
\right)$
\allowdisplaybreaks
\begin{align*}
&=
\Biggl\llbracket
\sigma^{\mathbf{PT}_{\boldsymbol{\mathcal{A}}^{(2)}}}
\left(\left(
Q_{j}
\circ^{0\mathbf{PT}_{\boldsymbol{\mathcal{A}}^{(2)}}}_{s_{j}}
P_{j}
\right)_{j\in\bb{\mathbf{s}}}
\right)
\Biggr\rrbracket_{s}
\tag{1}
\\&=
\Biggl\llbracket
\mathrm{CH}^{(2)}_{s}\left(
\mathrm{ip}^{(2,X)@}_{s}\left(
\sigma^{\mathbf{PT}_{\boldsymbol{\mathcal{A}}^{(2)}}}
\left(\left(
Q_{j}
\circ^{0\mathbf{PT}_{\boldsymbol{\mathcal{A}}^{(2)}}}_{s_{j}}
P_{j}
\right)_{j\in\bb{\mathbf{s}}}
\right)
\right)\right)
\Biggr\rrbracket_{s}
\tag{2}
\\&=
\Biggl\llbracket
\mathrm{CH}^{(2)}_{s}\left(
\sigma^{\mathbf{Pth}_{\boldsymbol{\mathcal{A}}^{(2)}}}
\left(\left(
\mathrm{ip}^{(2,X)@}_{s}\left(
Q_{j}
\right)
\circ^{0\mathbf{Pth}_{\boldsymbol{\mathcal{A}}^{(2)}}}_{s_{j}}
\mathrm{ip}^{(2,X)@}_{s}\left(
P_{j}
\right)
\right)_{j\in\bb{\mathbf{s}}}
\right)
\right)
\Biggr\rrbracket_{s}
\tag{3}
\\&=
\Biggl\llbracket
\mathrm{CH}^{(2)}_{s}\left(
\sigma^{\mathbf{Pth}_{\boldsymbol{\mathcal{A}}^{(2)}}}
\left(\left(
\mathrm{ip}^{(2,X)@}_{s}\left(
Q_{j}
\right)
\right)_{j\in\bb{\mathbf{s}}}
\right)
\circ^{0\mathbf{Pth}_{\boldsymbol{\mathcal{A}}^{(2)}}}_{s_{j}}
\right.
\\&\qquad\qquad\qquad\qquad\qquad\qquad\qquad\qquad
\left.
\sigma^{\mathbf{Pth}_{\boldsymbol{\mathcal{A}}^{(2)}}}
\left(\left(
\mathrm{ip}^{(2,X)@}_{s}\left(
P_{j}
\right)
\right)_{j\in\bb{\mathbf{s}}}
\right)
\right)
\Biggr\rrbracket_{s}
\tag{4}
\\&=
\Biggl\llbracket
\mathrm{CH}^{(2)}_{s}\left(
\mathrm{ip}^{(2,X)@}_{s}\left(
\sigma^{\mathbf{PT}_{\boldsymbol{\mathcal{A}}^{(2)}}}
\left(\left(
Q_{j}
\right)_{j\in\bb{\mathbf{s}}}\right)
\circ^{0\mathbf{PT}_{\boldsymbol{\mathcal{A}}^{(2)}}}
\sigma^{\mathbf{PT}_{\boldsymbol{\mathcal{A}}^{(2)}}}
\left(\left(
P_{j}
\right)_{j\in\bb{\mathbf{s}}}\right)
\right)\right)
\Biggr\rrbracket_{s}
\tag{5}
\\&=
\Biggl\llbracket
\sigma^{\mathbf{PT}_{\boldsymbol{\mathcal{A}}^{(2)}}}
\left(\left(
Q_{j}
\right)_{j\in\bb{\mathbf{s}}}\right)
\circ^{0\mathbf{PT}_{\boldsymbol{\mathcal{A}}^{(2)}}}
\sigma^{\mathbf{PT}_{\boldsymbol{\mathcal{A}}^{(2)}}}
\left(\left(
P_{j}
\right)_{j\in\bb{\mathbf{s}}}\right)
\Biggr\rrbracket_{s}
\tag{6}
\\&=
\sigma^{\llbracket\mathbf{PT}_{\boldsymbol{\mathcal{A}}^{(2)}}\rrbracket}
\left(\left(
\llbracket Q_{j}\rrbracket_{s_{j}}
\right)_{j\in\bb{\mathbf{s}}}\right)
\circ^{0\llbracket\mathbf{PT}_{\boldsymbol{\mathcal{A}}^{(2)}}\rrbracket}_{s}
\sigma^{\llbracket\mathbf{PT}_{\boldsymbol{\mathcal{A}}^{(2)}}\rrbracket}
\left(\left(
\llbracket P_{j}\rrbracket_{s_{j}}
\right)_{j\in\bb{\mathbf{s}}}\right).
\tag{7}
\end{align*}
\end{flushleft}

In the just stated chain of equalities, the first equality unravels the description of the operations of $0$-composition and $\sigma$ in the many sorted partial $\Sigma^{\boldsymbol{\mathcal{A}}^{(2)}}$-algebra $\llbracket\mathbf{PT}_{\boldsymbol{\mathcal{A}}^{(2)}}\rrbracket$ according to Proposition~\ref{PDPTQDCatAlg}; the second equality follows from Proposition~\ref{LDVWCong}; the third equality follows from the fact that $\mathrm{ip}^{(2,X)@}$ is a $\Sigma^{\boldsymbol{\mathcal{A}}^{(2)}}$-homomorphism according to Definition~\ref{DDIp}; the fourth equality follows from Proposition~\ref{PDVDCH} and Proposition~\ref{PDVVarA8}; the fifth equality follows from the fact that $\mathrm{ip}^{(2,X)@}$ is a $\Sigma^{\boldsymbol{\mathcal{A}}^{(2)}}$-homomorphism according to Definition~\ref{DDIp}; the sixth equality follows from Proposition~\ref{LDVWCong}; finally, the last equality recovers the description of the operations of $0$-composition and $\sigma$ in the many sorted partial $\Sigma^{\boldsymbol{\mathcal{A}}^{(2)}}$-algebra $\llbracket\mathbf{PT}_{\boldsymbol{\mathcal{A}}^{(2)}}\rrbracket$ according to Proposition~\ref{PDPTQDCatAlg}.

This completes the proof.
\end{proof}

\begin{proposition}\label{PDPTQVarB7} Let $(\mathbf{s},s)$ be an element in $S^{\star}\times S$, $\sigma$ an operation symbol in $\Sigma_{\mathbf{s},s}$ and 
$(\llbracket P_{j}\rrbracket_{s_{j}})_{j\in\bb{\mathbf{s}}}$ a family of second-order path term classes in  $\llbracket\mathrm{PT}_{\boldsymbol{\mathcal{A}}^{(2)}}\rrbracket_{\mathbf{s}}$, then the following equalities hold
\allowdisplaybreaks
\begin{align*}
\sigma^{\llbracket\mathbf{PT}_{\boldsymbol{\mathcal{A}}^{(2)}}\rrbracket}
\left(\left(
\mathrm{sc}^{1\llbracket\mathbf{PT}_{\boldsymbol{\mathcal{A}}^{(2)}}\rrbracket}_{s_{j}}
\left(
\llbracket P_{j}\rrbracket_{s_{j}}
\right)\right)_{j\in\bb{\mathbf{s}}}
\right)
&=
\mathrm{sc}^{1\llbracket\mathbf{PT}_{\boldsymbol{\mathcal{A}}^{(2)}}\rrbracket}_{s}\left(
\sigma^{\llbracket\mathbf{PT}_{\boldsymbol{\mathcal{A}}^{(2)}}\rrbracket}
\left(\left(
\llbracket P_{j}\rrbracket_{s_{j}}
\right)_{j\in\bb{\mathbf{s}}}
\right)
\right);
\\
\sigma^{\llbracket\mathbf{PT}_{\boldsymbol{\mathcal{A}}^{(2)}}\rrbracket}
\left(\left(
\mathrm{tg}^{1\llbracket\mathbf{PT}_{\boldsymbol{\mathcal{A}}^{(2)}}\rrbracket}_{s_{j}}
\left(
\llbracket P_{j}\rrbracket_{s_{j}}
\right)\right)_{j\in\bb{\mathbf{s}}}
\right)
&=
\mathrm{tg}^{1\llbracket\mathbf{PT}_{\boldsymbol{\mathcal{A}}^{(2)}}\rrbracket}_{s}\left(
\sigma^{\llbracket\mathbf{PT}_{\boldsymbol{\mathcal{A}}^{(2)}}\rrbracket}
\left(\left(
\llbracket P_{j}\rrbracket_{s_{j}}
\right)_{j\in\bb{\mathbf{s}}}
\right)
\right).
\end{align*}
\end{proposition}
\begin{proof}
We will only present the proof for the first equality. The other one is handled in a similar manner.

Note that the following chain of equalities holds
\begin{flushleft}
$\sigma^{\llbracket\mathbf{PT}_{\boldsymbol{\mathcal{A}}^{(2)}}\rrbracket}
\left(\left(
\mathrm{sc}^{1\llbracket\mathbf{PT}_{\boldsymbol{\mathcal{A}}^{(2)}}\rrbracket}_{s_{j}}
\left(
\llbracket P_{j}\rrbracket_{s_{j}}
\right)\right)_{j\in\bb{\mathbf{s}}}
\right)$
\allowdisplaybreaks
\begin{align*}
&=
\Biggl\llbracket
\sigma^{\mathbf{PT}_{\boldsymbol{\mathcal{A}}^{(2)}}}
\left(\left(
\mathrm{sc}^{1\mathbf{PT}_{\boldsymbol{\mathcal{A}}^{(2)}}}_{s_{j}}
\left(
P_{j}
\right)
\right)_{j\in\bb{\mathbf{s}}}
\right)
\Biggr\rrbracket_{s}
\tag{1}
\\
&=
\Biggl\llbracket
\mathrm{CH}^{(2)}_{s}\left(
\mathrm{ip}^{(2,X)@}_{s}\left(
\sigma^{\mathbf{PT}_{\boldsymbol{\mathcal{A}}^{(2)}}}
\left(\left(
\mathrm{sc}^{1\mathbf{PT}_{\boldsymbol{\mathcal{A}}^{(2)}}}_{s_{j}}
\left(
P_{j}
\right)
\right)_{j\in\bb{\mathbf{s}}}
\right)
\right)
\right)
\Biggr\rrbracket_{s}
\tag{2}
\\
&=
\Biggl\llbracket
\mathrm{CH}^{(2)}_{s}\left(
\sigma^{\mathbf{Pth}_{\boldsymbol{\mathcal{A}}^{(2)}}}
\left(\left(
\mathrm{sc}^{1\mathbf{Pth}_{\boldsymbol{\mathcal{A}}^{(2)}}}_{s_{j}}
\left(
\mathrm{ip}^{(2,X)@}_{s_{j}}\left(
P_{j}
\right)
\right)
\right)_{j\in\bb{\mathbf{s}}}
\right)
\right)
\Biggr\rrbracket_{s}
\tag{3}
\\
&=
\Biggl\llbracket
\mathrm{CH}^{(2)}_{s}\left(
\mathrm{sc}^{1\mathbf{Pth}_{\boldsymbol{\mathcal{A}}^{(2)}}}_{s}
\left(
\sigma^{\mathbf{Pth}_{\boldsymbol{\mathcal{A}}^{(2)}}}
\left(\left(
\mathrm{ip}^{(2,X)@}_{s_{j}}\left(
P_{j}
\right)
\right)_{j\in\bb{\mathbf{s}}}
\right)
\right)
\right)
\Biggr\rrbracket_{s}
\tag{4}
\\
&=
\Biggl\llbracket
\mathrm{CH}^{(2)}_{s}\left(
\mathrm{ip}^{(2,X)@}_{s}\left(
\mathrm{sc}^{1\mathbf{PT}_{\boldsymbol{\mathcal{A}}^{(2)}}}_{s}
\left(
\sigma^{\mathbf{PT}_{\boldsymbol{\mathcal{A}}^{(2)}}}
\left(\left(
P_{j}
\right)_{j\in\bb{\mathbf{s}}}
\right)
\right)
\right)
\right)
\Biggr\rrbracket_{s}
\tag{5}
\\
&=
\Biggl\llbracket
\mathrm{sc}^{1\mathbf{PT}_{\boldsymbol{\mathcal{A}}^{(2)}}}_{s}
\left(
\sigma^{\mathbf{PT}_{\boldsymbol{\mathcal{A}}^{(2)}}}
\left(\left(
P_{j}
\right)_{j\in\bb{\mathbf{s}}}
\right)
\right)
\Biggr\rrbracket_{s}
\tag{6}
\\&=
\mathrm{sc}^{1\llbracket\mathbf{PT}_{\boldsymbol{\mathcal{A}}^{(2)}}\rrbracket}_{s}\left(
\sigma^{\llbracket\mathbf{PT}_{\boldsymbol{\mathcal{A}}^{(2)}}\rrbracket}
\left(\left(
\llbracket P_{j}\rrbracket_{s_{j}}
\right)_{j\in\bb{\mathbf{s}}}
\right)
\right).
\tag{7}
\end{align*}
\end{flushleft}

In the just stated chain of equalities, the first equality unravels the description of the operations of $1$-source and $\sigma$ in the many sorted partial $\Sigma^{\boldsymbol{\mathcal{A}}^{(2)}}$-algebra $\llbracket\mathbf{PT}_{\boldsymbol{\mathcal{A}}^{(2)}}\rrbracket$ according to Proposition~\ref{PDPTQDCatAlg}; the second equality follows from Proposition~\ref{LDVWCong}; the third equality follows from the fact that $\mathrm{ip}^{(2,X)@}$ is a $\Sigma^{\boldsymbol{\mathcal{A}}^{(2)}}$-homomorphism according to Definition~\ref{DDIp}; the fourth equality follows from Proposition~\ref{PDVDCH} and Proposition~\ref{PDVVarB7}; the fifth equality follows from the fact that $\mathrm{ip}^{(2,X)@}$ is a $\Sigma^{\boldsymbol{\mathcal{A}}^{(2)}}$-homomorphism according to Definition~\ref{DDIp}; the sixth equality follows from Proposition~\ref{LDVWCong}; finally, the last equality recovers the description of the operations of $1$-source and $\sigma$ in the many sorted partial $\Sigma^{\boldsymbol{\mathcal{A}}^{(2)}}$-algebra $\llbracket\mathbf{PT}_{\boldsymbol{\mathcal{A}}^{(2)}}\rrbracket$ according to Proposition~\ref{PDPTQDCatAlg}.

This completes the proof.
\end{proof}

\begin{proposition}\label{PDPTQVarB8} Let $(\mathbf{s},s)$ be an element in $S^{\star}\times S$, $\sigma$ and operation symbol in $\Sigma_{\mathbf{s},s}$ and $(\llbracket Q_{j}\rrbracket_{s_{j}})_{j\in\bb{\mathbf{s}}}$, $(\llbracket P_{j}\rrbracket_{s_{j}})_{j\in\bb{\mathbf{s}}}$ be two families of second-order path term classes in  $\llbracket\mathrm{PT}_{\boldsymbol{\mathcal{A}}^{(2)}}\rrbracket_{\mathbf{s}}$ satisfying that, for every $j\in\bb{\mathbf{s}}$,
\[\mathrm{sc}^{1\llbracket\mathbf{PT}_{\boldsymbol{\mathcal{A}}^{(2)}}\rrbracket}_{s_{j}}\left(
\llbracket Q_{j}\rrbracket_{s_{j}}
\right)=
\mathrm{tg}^{1\llbracket\mathbf{PT}_{\boldsymbol{\mathcal{A}}^{(2)}}\rrbracket}_{s_{j}}\left(
\llbracket P_{j}\rrbracket_{s_{j}}
\right),
\]
then the following equality holds
\allowdisplaybreaks
\begin{multline*}
\sigma^{\llbracket\mathbf{PT}_{\boldsymbol{\mathcal{A}}^{(2)}}\rrbracket}
\left(\left(
\llbracket Q_{j}\rrbracket_{s_{j}}
\circ^{1\llbracket\mathbf{PT}_{\boldsymbol{\mathcal{A}}^{(2)}}\rrbracket}_{s_{j}}
\llbracket P_{j}\rrbracket_{s_{j}}
\right)_{j\in\bb{\mathbf{s}}}
\right)
\\=
\sigma^{\llbracket\mathbf{PT}_{\boldsymbol{\mathcal{A}}^{(2)}}\rrbracket}
\left(\left(
\llbracket Q_{j}\rrbracket_{s_{j}}
\right)_{j\in\bb{\mathbf{s}}}\right)
\circ^{1\llbracket\mathbf{PT}_{\boldsymbol{\mathcal{A}}^{(2)}}\rrbracket}_{s}
\sigma^{\llbracket\mathbf{PT}_{\boldsymbol{\mathcal{A}}^{(2)}}\rrbracket}
\left(\left(
\llbracket P_{j}\rrbracket_{s_{j}}
\right)_{j\in\bb{\mathbf{s}}}\right).
\end{multline*}
\end{proposition}
\begin{proof}
Note that the following chain of equalities holds
\begin{flushleft}
$\sigma^{\llbracket\mathbf{PT}_{\boldsymbol{\mathcal{A}}^{(2)}}\rrbracket}
\left(\left(
\llbracket Q_{j}\rrbracket_{s_{j}}
\circ^{1\llbracket\mathbf{PT}_{\boldsymbol{\mathcal{A}}^{(2)}}\rrbracket}_{s_{j}}
\llbracket P_{j}\rrbracket_{s_{j}}
\right)_{j\in\bb{\mathbf{s}}}
\right)$
\allowdisplaybreaks
\begin{align*}
&=
\Biggl\llbracket
\sigma^{\mathbf{PT}_{\boldsymbol{\mathcal{A}}^{(2)}}}
\left(\left(
Q_{j}
\circ^{1\mathbf{PT}_{\boldsymbol{\mathcal{A}}^{(2)}}}_{s_{j}}
P_{j}
\right)_{j\in\bb{\mathbf{s}}}
\right)
\Biggr\rrbracket_{s}
\tag{1}
\\&=
\Biggl\llbracket
\mathrm{CH}^{(2)}_{s}\left(
\mathrm{ip}^{(2,X)@}_{s}\left(
\sigma^{\mathbf{PT}_{\boldsymbol{\mathcal{A}}^{(2)}}}
\left(\left(
Q_{j}
\circ^{1\mathbf{PT}_{\boldsymbol{\mathcal{A}}^{(2)}}}_{s_{j}}
P_{j}
\right)_{j\in\bb{\mathbf{s}}}
\right)
\right)\right)
\Biggr\rrbracket_{s}
\tag{2}
\\&=
\Biggl\llbracket
\mathrm{CH}^{(2)}_{s}\left(
\sigma^{\mathbf{Pth}_{\boldsymbol{\mathcal{A}}^{(2)}}}
\left(\left(
\mathrm{ip}^{(2,X)@}_{s}\left(
Q_{j}
\right)
\circ^{1\mathbf{Pth}_{\boldsymbol{\mathcal{A}}^{(2)}}}_{s_{j}}
\mathrm{ip}^{(2,X)@}_{s}\left(
P_{j}
\right)
\right)_{j\in\bb{\mathbf{s}}}
\right)
\right)
\Biggr\rrbracket_{s}
\tag{3}
\\&=
\Biggl\llbracket
\mathrm{CH}^{(2)}_{s}\left(
\sigma^{\mathbf{Pth}_{\boldsymbol{\mathcal{A}}^{(2)}}}
\left(\left(
\mathrm{ip}^{(2,X)@}_{s}\left(
Q_{j}
\right)
\right)_{j\in\bb{\mathbf{s}}}
\right)
\circ^{1\mathbf{Pth}_{\boldsymbol{\mathcal{A}}^{(2)}}}_{s_{j}}
\right.
\\&\qquad\qquad\qquad\qquad\qquad\qquad\qquad\qquad
\left.
\sigma^{\mathbf{Pth}_{\boldsymbol{\mathcal{A}}^{(2)}}}
\left(\left(
\mathrm{ip}^{(2,X)@}_{s}\left(
P_{j}
\right)
\right)_{j\in\bb{\mathbf{s}}}
\right)
\right)
\Biggr\rrbracket_{s}
\tag{4}
\\&=
\Biggl\llbracket
\mathrm{CH}^{(2)}_{s}\left(
\mathrm{ip}^{(2,X)@}_{s}\left(
\sigma^{\mathbf{PT}_{\boldsymbol{\mathcal{A}}^{(2)}}}
\left(\left(
Q_{j}
\right)_{j\in\bb{\mathbf{s}}}\right)
\circ^{1\mathbf{PT}_{\boldsymbol{\mathcal{A}}^{(2)}}}
\sigma^{\mathbf{PT}_{\boldsymbol{\mathcal{A}}^{(2)}}}
\left(\left(
P_{j}
\right)_{j\in\bb{\mathbf{s}}}\right)
\right)\right)
\Biggr\rrbracket_{s}
\tag{5}
\\&=
\Biggl\llbracket
\sigma^{\mathbf{PT}_{\boldsymbol{\mathcal{A}}^{(2)}}}
\left(\left(
Q_{j}
\right)_{j\in\bb{\mathbf{s}}}\right)
\circ^{1\mathbf{PT}_{\boldsymbol{\mathcal{A}}^{(2)}}}
\sigma^{\mathbf{PT}_{\boldsymbol{\mathcal{A}}^{(2)}}}
\left(\left(
P_{j}
\right)_{j\in\bb{\mathbf{s}}}\right)
\Biggr\rrbracket_{s}
\tag{6}
\\&=
\sigma^{\llbracket\mathbf{PT}_{\boldsymbol{\mathcal{A}}^{(2)}}\rrbracket}
\left(\left(
\llbracket Q_{j}\rrbracket_{s_{j}}
\right)_{j\in\bb{\mathbf{s}}}\right)
\circ^{1\llbracket\mathbf{PT}_{\boldsymbol{\mathcal{A}}^{(2)}}\rrbracket}_{s}
\sigma^{\llbracket\mathbf{PT}_{\boldsymbol{\mathcal{A}}^{(2)}}\rrbracket}
\left(\left(
\llbracket P_{j}\rrbracket_{s_{j}}
\right)_{j\in\bb{\mathbf{s}}}\right).
\tag{7}
\end{align*}
\end{flushleft}

In the just stated chain of equalities, the first equality unravels the description of the operations of $1$-composition and $\sigma$ in the many sorted partial $\Sigma^{\boldsymbol{\mathcal{A}}^{(2)}}$-algebra $\llbracket\mathbf{PT}_{\boldsymbol{\mathcal{A}}^{(2)}}\rrbracket$ according to Proposition~\ref{PDPTQDCatAlg}; the second equality follows from Proposition~\ref{LDVWCong}; the third equality follows from the fact that $\mathrm{ip}^{(2,X)@}$ is a $\Sigma^{\boldsymbol{\mathcal{A}}^{(2)}}$-homomorphism according to Definition~\ref{DDIp}; the fourth equality follows from Proposition~\ref{PDVDCH} and Proposition~\ref{PDVVarB8}; the fifth equality follows from the fact that $\mathrm{ip}^{(2,X)@}$ is a $\Sigma^{\boldsymbol{\mathcal{A}}^{(2)}}$-homomorphism according to Definition~\ref{DDIp}; the sixth equality follows from Proposition~\ref{LDVWCong}; finally, the last equality recovers the description of the operations of $1$-composition and $\sigma$ in the many sorted partial $\Sigma^{\boldsymbol{\mathcal{A}}^{(2)}}$-algebra $\llbracket\mathbf{PT}_{\boldsymbol{\mathcal{A}}^{(2)}}\rrbracket$ according to Proposition~\ref{PDPTQDCatAlg}.

This completes the proof.
\end{proof}

\begin{restatable}{proposition}{PDPTQDCtyAlg}
\label{PDPTQDCtyAlg} $\llbracket\mathsf{PT}_{\boldsymbol{\mathcal{A}}^{(2)}}\rrbracket$ is a $2$-categorial $\Sigma$-algebra.
\end{restatable}
\begin{proof}
That  $\llbracket\mathsf{PT}_{\boldsymbol{\mathcal{A}}^{(2)}}\rrbracket$ is an $S$-sorted $2$-category was already proven in Proposition~\ref{PDPTQDCat}. We are only left to prove that, for every $(\mathbf{s},s)$ in $S^{\star}\times S$, $\sigma^{\llbracket\mathbf{PT}_{\boldsymbol{\mathcal{A}}^{(2)}}\rrbracket}$ is a $2$-functor from $\llbracket\mathsf{PT}_{\boldsymbol{\mathcal{A}}^{(2)}}\rrbracket_{\mathbf{s}}$ to $\llbracket\mathsf{PT}_{\boldsymbol{\mathcal{A}}^{(2)}}\rrbracket_{s}$. This is a fact that follows from Propositions~\ref{PDPTQVarA7},~\ref{PDPTQVarA8},~\ref{PDPTQVarB7} and~\ref{PDPTQVarB8}.
\end{proof}

\section{
\texorpdfstring
{An Artinian preorder on $\coprod \llbracket\mathrm{PT}_{\boldsymbol{\mathcal{A}}^{(2)}}\rrbracket$}
{An Artinian preorder on the quotient of second-order path terms}
}
In this section we define on $\coprod \llbracket\mathrm{PT}_{\boldsymbol{\mathcal{A}}^{(2)}}\rrbracket$, the coproduct of the many-sorted set of second-order path term classes---formed by all labelled second-order path term classes $(\llbracket P\rrbracket_{s},s)$ with $s\in S$ and $\llbracket P\rrbracket_{s}\in\llbracket\mathrm{PT}_{\boldsymbol{\mathcal{A}}^{(2)}}\rrbracket_{s}$---, an Artinian preorder. 

\begin{restatable}{definition}{DDPTQOrd}
\label{DDPTQOrd} 
\index{partial preorder!second-order!$\leq_{\llbracket\mathbf{PT}_{\boldsymbol{\mathcal{A}}^{(2)}}\rrbracket}$}
We let $\leq_{\llbracket\mathbf{PT}_{\boldsymbol{\mathcal{A}}^{(2)}}\rrbracket}$ stand for the relation defined on $\coprod \llbracket\mathrm{PT}_{\boldsymbol{\mathcal{A}}^{(2)}}\rrbracket$ which contains the ordered pairs
$((\llbracket Q\rrbracket_{t},t),(\llbracket P\rrbracket_{s},s))$ in $(\coprod \llbracket\mathrm{PT}_{\boldsymbol{\mathcal{A}}^{(2)}}\rrbracket)^{2}$ such that
$$
\left(
\Bigl\llbracket
\mathrm{ip}^{(2,X)@}_{t}\left(
Q
\right)
\Bigr\rrbracket_{t}
,t\right)
\leq_{\llbracket\mathbf{Pth}_{\boldsymbol{\mathcal{A}}^{(2)}}\rrbracket}
\left(
\Bigl\llbracket
\mathrm{ip}^{(2,X)@}_{s}\left(
P
\right)\Bigr\rrbracket_{s},s\right).
$$

Thus, $(\llbracket Q\rrbracket_{t},t)$ $\leq_{\llbracket\mathbf{PT}_{\boldsymbol{\mathcal{A}}^{(2)}}\rrbracket}$-precedes $(\llbracket P\rrbracket_{s},s)$ if the $\llbracket \cdot \rrbracket$-class of the second-order path obtained from applying the $\mathrm{ip}^{(2,X)@}$ mapping at $Q$ $\leq_{\llbracket\mathbf{Pth}_{\boldsymbol{\mathcal{A}}^{(2)}}\rrbracket}$-precedes the respective second-order path obtained from applying the  $\mathrm{ip}^{(2,X)@}$ mapping at $P$.
\end{restatable}

\begin{remark}
Note that the relation $\leq_{\llbracket\mathbf{PT}_{\boldsymbol{\mathcal{A}}^{(2)}}\rrbracket}$ on $\coprod \llbracket\mathrm{PT}_{\boldsymbol{\mathcal{A}}^{(2)}}\rrbracket$ is well-defined since, for every sort $s\in S$ and every pair of second-order path terms $P,P'\in\mathrm{PT}_{\boldsymbol{\mathcal{A}}^{(2)},s}$ satisfying that $(P,P')\in\Theta^{\llbracket 2 \rrbracket}_{s}$, we have, by Lemma~\ref{LDVCong}, that 
$$
\Bigl\llbracket
\mathrm{ip}^{(2,X)@}_{s}\left(
P
\right)
\Bigr\rrbracket_{s}=
\Bigl\llbracket
\mathrm{ip}^{(2,X)@}_{s}\left(
P'
\right)
\Bigr\rrbracket_{s}.
$$
\end{remark}

We next prove that the just defined relation is an Artinian preorder on $\coprod \llbracket\mathrm{PT}_{\boldsymbol{\mathcal{A}}^{(2)}}\rrbracket$.

\begin{restatable}{proposition}{PDPTQOrdArt}
\label{PDPTQOrdArt} $(\coprod\llbracket\mathrm{PT}_{\boldsymbol{\mathcal{A}}^{(2)}}\rrbracket, \leq_{\llbracket\mathbf{PT}_{\boldsymbol{\mathcal{A}}^{(2)}}\rrbracket})$ is a preordered set. Moreover, in this preordered set there is not any strictly decreasing $\omega_{0}$-chain, i.e., $(\coprod\llbracket\mathrm{PT}_{\boldsymbol{\mathcal{A}}^{(2)}}\rrbracket, \leq_{\llbracket\mathbf{PT}_{\boldsymbol{\mathcal{A}}^{(2)}}\rrbracket})$ is an Artinian partially preordered set.
\end{restatable}
\begin{proof}
Note that for every sort $s\in S$ and every second-order path term $P\in\mathrm{PT}_{\boldsymbol{\mathcal{A}}^{(2)},s}$, the class $\llbracket \mathrm{ip}^{(2,X)@}_{s}(P) \rrbracket_{s}$ belongs to $\llbracket \mathrm{Pth}_{\boldsymbol{\mathcal{A}}^{(2)}} \rrbracket_{s}$. The statement follows from the fact that $(\coprod\llbracket \mathrm{Pth}_{\boldsymbol{\mathcal{A}}^{(2)}} \rrbracket ,\leq_{\llbracket \mathrm{Pth}_{\boldsymbol{\mathcal{A}}^{(2)}} \rrbracket})$ is an Artinian preordered set according to Proposition~\ref{PDVOrdArt}.
\end{proof}			
\chapter{Second-order isomorphisms}\label{S2N}

In this chapter we prove that the many-sorted partial $\Sigma^{\boldsymbol{\mathcal{A}}^{(2)}}$-algebras of second-order path classes $\llbracket \mathbf{Pth}_{\boldsymbol{\mathcal{A}}^{(2)}}\rrbracket$ and second-order path term classes $\llbracket \mathbf{PT}_{\boldsymbol{\mathcal{A}}^{(2)}}\rrbracket$ are isomorphic. Indeed the mappings $\mathrm{ip}^{(\llbracket 2 \rrbracket,X)@}$ and $\mathrm{CH}^{\llbracket 2 \rrbracket}$ introduced in the previous chapter form a pair of inverse $\Sigma^{\boldsymbol{\mathcal{A}}^{(2)}}$-isomorphisms. Moreover, it is also shown that the many-sorted $2$-categorial $\Sigma$-algebras $\llbracket \mathsf{Pth}_{\boldsymbol{\mathcal{A}}^{(2)}}\rrbracket$ and $\llbracket\mathsf{PT}_{\boldsymbol{\mathcal{A}}^{(2)}}\rrbracket$ are also isomorphic, since the mappings $\mathrm{ip}^{(\llbracket 2 \rrbracket,X)@}$ and $\mathrm{CH}^{\llbracket 2 \rrbracket}$ form also a pair of inverse functors, i.e., of $2$-categorial $\Sigma$-isomorphisms. Moreover, the mappings $\mathrm{ip}^{(\llbracket 2 \rrbracket,X)@}$ and $\mathrm{CH}^{\llbracket 2 \rrbracket}$ form also a pair of inverse order-preserving mappings which, in turn, implies that the Artinian partial preorders $\leq_{\llbracket \mathbf{PT}_{\boldsymbol{\mathcal{A}}^{(2)}}\rrbracket}$ and $\leq_{\llbracket \mathbf{Pth}_{\boldsymbol{\mathcal{A}}^{(2)}}\rrbracket}$ are isomorphic.


\section{A second-order algebraic isomorphism}

We begin by showing that the many-sorted partial $\Sigma^{\boldsymbol{\mathcal{A}}^{(2)}}$-algebras $\llbracket\mathbf{Pth}_{\boldsymbol{\mathcal{A}}^{(2)}}\rrbracket$ of second-order path classes, introduced in Proposition~\ref{PDVDCatAlg}, and 
$\llbracket\mathbf{PT}_{\boldsymbol{\mathcal{A}}^{(2)}}\rrbracket$ of second-order path term classes, introduced in Proposition~\ref{PDPTQDCatAlg}, are isomorphic.

\begin{restatable}{theorem}{TDIso}
\label{TDIso} 
The many-sorted partial $\Sigma^{\boldsymbol{\mathcal{A}}^{(2)}}$-algebras 
$\llbracket\mathbf{Pth}_{\boldsymbol{\mathcal{A}}^{(2)}}\rrbracket$ and $\llbracket\mathbf{PT}_{\boldsymbol{\mathcal{A}}^{(2)}}\rrbracket$ are isomorphic.
\end{restatable}
\begin{proof} We claim that the mappings
\begin{center}
\begin{tikzpicture}
[ACliment/.style={-{To [angle'=45, length=5.75pt, width=4pt, round]}}
, scale=0.8, 
AClimentD/.style={double equal sign distance,
-implies
}
]
\node[] (PK) at (-2.9,-3) [] {$\llbracket\mathbf{Pth}_{\boldsymbol{\mathcal{A}}^{(2)}}\rrbracket$};
\node[] (PTK) at (2.9,-3) [] {$\llbracket\mathbf{PT}_{\boldsymbol{\mathcal{A}}^{(2)}}\rrbracket$};
 
 \draw[shorten >=0.2cm,shorten <=.2cm, ACliment]
 ($(-2,-3)+.09*(0,3)$) to node  [above]  {$ \mathrm{CH}^{\llbracket 2\rrbracket}$}  ($(2,-3)+.09*(0,3)$);
 \draw[shorten >=0.2cm,shorten <=.2cm, ACliment]
 ($(2,-3)-.09*(0,3)$) to node  [below]  {$\mathrm{ip}^{(\llbracket 2 \rrbracket,X)@}$}  ($(-2,-3)-.09*(0,3)$);
\end{tikzpicture}
\end{center}
define a pair of $\Sigma^{\boldsymbol{\mathcal{A}}^{(2)}}$-isomorphisms satisfying that 
\allowdisplaybreaks
\begin{align*}
\mathrm{ip}^{(\llbracket 2 \rrbracket,X)@}
\circ
\mathrm{CH}^{\llbracket 2\rrbracket}
&=
\mathrm{id}^{\llbracket\mathbf{Pth}_{\boldsymbol{\mathcal{A}}^{(2)}}\rrbracket};
&
\mathrm{CH}^{\llbracket 2\rrbracket}
\circ
\mathrm{ip}^{(\llbracket 2 \rrbracket,X)@}
&=
\mathrm{id}^{\llbracket\mathbf{PT}_{\boldsymbol{\mathcal{A}}^{(2)}}\rrbracket}.
\end{align*}

We prove some intermediate claims that will allows us to conclude the aforementioned statement.

\begin{restatable}{claim}{CDIsoCH}
\label{CDIsoCH}
The mapping
$\mathrm{CH}^{\llbracket 2\rrbracket}$
is a $\Sigma^{\boldsymbol{\mathcal{A}}^{(2)}}$-isomorphism.
\end{restatable}

Let us recall that $\mathrm{CH}^{\llbracket 2\rrbracket}$, introduced in Definition~\ref{DDPTQDCH}, is a many-sorted mapping of the form
$$
\mathrm{CH}^{\llbracket 2\rrbracket}
\colon
\llbracket \mathrm{Pth}_{\boldsymbol{\mathcal{A}}^{(2)}}\rrbracket
\mor
\llbracket\mathrm{PT}_{\boldsymbol{\mathcal{A}}^{(2)}}\rrbracket.
$$

Following Definition~\ref{DDPTQDCH}, the many-sorted mapping $\mathrm{CH}^{\llbracket 2\rrbracket}$ is given by
\[
\mathrm{CH}^{\llbracket 2\rrbracket}=\mathrm{pr}^{\Theta^{\llbracket 2\rrbracket}}
\circ \mathrm{CH}^{(2)\mathrm{m}}
=
\mathrm{p}^{(\Theta^{\llbracket 2\rrbracket},\Theta^{[2]})}
\circ
\mathrm{pr}^{\Theta^{[2]}}
\circ
\mathrm{CH}^{(2)}
\circ
\mathrm{pr}^{\mathrm{Ker}(\mathrm{CH}^{(2)})}
.
\]

Following Proposition~\ref{PDCHDCatAlg}, $\mathrm{pr}^{\mathrm{Ker}(\mathrm{CH}^{(2)})}$ is a $\Sigma^{\boldsymbol{\mathcal{A}}^{(2)}}$-homomorphism. Following Proposition~\ref{PDThetaDCH}, $\mathrm{pr}^{\Theta^{[2]}}\circ \mathrm{CH}^{(2)}$ is a $\Sigma^{\boldsymbol{\mathcal{A}}}$-homomorphism. Moreover, the mapping $\mathrm{p}^{(\Theta^{\llbracket 2\rrbracket},\Theta^{[2]})}$ is a $\Sigma^{\boldsymbol{\mathcal{A}}^{(2)}}$-homomorphism. It follows that  $\mathrm{CH}^{\llbracket 2\rrbracket}$, being a composition of $\Sigma^{\boldsymbol{\mathcal{A}}^{(2)}}$-homomorphisms, is itself a $\Sigma^{\boldsymbol{\mathcal{A}}^{(2)}}$-homomorphism.

This finishes the proof of Claim~\ref{CDIsoCH}.

Now, we  turn ourselves to study the case of the many-sorted mapping $\mathrm{ip}^{(\llbracket 2 \rrbracket,X)@}$.

\begin{restatable}{claim}{CDIsoIpfc}
\label{CDIsoIpfc} The mapping $\mathrm{ip}^{(\llbracket 2 \rrbracket,X)@}$ is a $\Sigma^{\boldsymbol{\mathcal{A}}^{(2)}}$-isomorphism.
\end{restatable}

Let us recall that $\mathrm{ip}^{(\llbracket 2 \rrbracket,X)@}$, introduced in Definition~\ref{DDPTQIp}, is a many-sorted mapping of the form
$$
\mathrm{ip}^{(\llbracket 2 \rrbracket,X)@}
\colon
\llbracket\mathrm{PT}_{\boldsymbol{\mathcal{A}}^{(2)}}\rrbracket
\mor
\llbracket \mathrm{Pth}_{\boldsymbol{\mathcal{A}}^{(2)}}\rrbracket.
$$

Following Definition~\ref{DDPTQIp}, the many-sorted mapping $
\mathrm{ip}^{(\llbracket 2 \rrbracket,X)@}$ is given by
\[
\mathrm{ip}^{(\llbracket 2 \rrbracket,X)@}=
\left(\mathrm{pr}^{\llbracket \cdot \rrbracket}\circ\mathrm{ip}^{(2,X)@}\right)^{\mathrm{m}}\circ\mathrm{p}^{
\mathrm{pr}^{\llbracket \cdot \rrbracket}\circ\mathrm{ip}^{(2,X)@},
\Theta^{\llbracket 2 \rrbracket}
}.
\]

Following Definition~\ref{DDPTQIp}, this mapping satisfies that 
\[
\mathrm{pr}^{\llbracket \cdot \rrbracket}\circ \mathrm{ip}^{(2,X)@}=
\mathrm{p}^{
\mathrm{pr}^{\llbracket \cdot \rrbracket}\circ\mathrm{ip}^{(2,X)@},
\Theta^{\llbracket 2 \rrbracket}
}
\circ \mathrm{pr}^{\Theta^{\llbracket 2\rrbracket}}.
\]

Note that $\mathrm{ip}^{(2,X)@}$ is a $\Sigma^{\boldsymbol{\mathcal{A}}^{(2)}}$-homomorphism according to Definition~\ref{DDIp} and $\mathrm{pr}^{\llbracket \cdot \rrbracket}$ is a $\Sigma^{\boldsymbol{\mathcal{A}}^{(2)}}$-homomorphism according to Proposition~\ref{PDVDCatAlg}. It follows that  $\mathrm{ip}^{(\llbracket 2 \rrbracket,X)@}$, being a composition of $\Sigma^{\boldsymbol{\mathcal{A}}^{(2)}}$-homomorphisms, is itself a $\Sigma^{\boldsymbol{\mathcal{A}}^{(2)}}$-homomorphism.

This finishes the proof of Claim~\ref{CDIsoIpfc}.

We end by showing that this pair of isomorphisms are reciprocally inverse.

\begin{restatable}{claim}{CDIso}
\label{CDIso}
$\mathrm{CH}^{\llbracket 2\rrbracket}$ and $\mathrm{ip}^{(\llbracket 2 \rrbracket,X)@}$ form a pair of inverse $\Sigma^{\boldsymbol{\mathcal{A}}^{(2)}}$-isomorphisms.
\end{restatable}

Let $s$ be a sort in $S$ and $\llbracket \mathfrak{P}^{(2)}\rrbracket_{s}$ a second-order path class in $\llbracket \mathrm{Pth}_{\boldsymbol{\mathcal{A}}^{(2)}}\rrbracket_{s}$.

The following chain of equalities holds
\begin{align*}
\mathrm{ip}^{(\llbracket 2 \rrbracket,X)@}_{s}\left(
\mathrm{CH}^{\llbracket 2\rrbracket}_{s}\left(
\Bigl\llbracket
\mathfrak{P}^{(2)}
\Bigr\rrbracket_{s}
\right)\right)
&=
\biggl\llbracket
\mathrm{ip}^{(2,X)@}_{s}\left(
\mathrm{CH}^{(2)}_{s}\left(
\mathfrak{P}^{(2)}
\right)\right)
\biggr\rrbracket_{s}
\tag{1}
\\&=
\Bigl\llbracket
\mathfrak{P}^{(2)}
\Bigr\rrbracket_{s}.
\tag{2}
\end{align*}

The first equality unravels the definition of the mappings $\mathrm{ip}^{(\llbracket 2 \rrbracket,X)@}$ and $\mathrm{CH}^{\llbracket 2\rrbracket}$ according to, respectively, Definitions~\ref{DDPTQIp} and~\ref{DDPTQDCH}; finally, the second equality follows from Proposition~\ref{PDIpDCH}.

On the other hand, let $s$ be a sort in $S$ and $\llbracket P\rrbracket_{s}$ a second-order path term class in $\llbracket\mathbf{PT}_{\boldsymbol{\mathcal{A}}^{(2)}}\rrbracket_{s}$.

The following chain of equalities holds
\begin{align*}
\mathrm{CH}^{\llbracket 2\rrbracket}_{s}\left(
\mathrm{ip}^{(\llbracket 2 \rrbracket,X)@}_{s}\left(
\llbracket
P
\rrbracket_{s}
\right)\right)
&=
\Bigl\llbracket
\mathrm{CH}^{(2)}_{s}\left(
\mathrm{ip}^{2,X,@}_{s}\left(
P
\right)\right)
\Bigr\rrbracket_{s}
\tag{1}
\\&=
\llbracket
P
\rrbracket_{s}.
\tag{2}
\end{align*}

The first equality unravels  the definition of the mappings $\mathrm{ip}^{(\llbracket 2 \rrbracket,X)@}$ and $\mathrm{CH}^{\llbracket 2\rrbracket}$ according to, respectively, Definitions~\ref{DDPTQIp} and~\ref{DDPTQDCH}; finally, the second equality follows from Lemma~\ref{LDVWCong}.

This completes the proof of Claim~\ref{CDIso}.

This completes the proof of Theorem~\ref{TDIso}.
\end{proof}

\section{A  second-order categorial  isomorphism}
In this section we prove that the $S$-sorted categories $\llbracket \mathsf{Pth}_{\boldsymbol{\mathcal{A}}^{(2)}}\rrbracket$, introduced in Definition~\ref{PDVCat}, and $\llbracket\mathsf{PT}_{\boldsymbol{\mathcal{A}}^{(2)}}\rrbracket$, introduced in Definition~\ref{DDPTQDCat}, are isomorphic.

\begin{restatable}{theorem}{TDIsoCat}
\label{TDIsoCat} The many-sorted $2$-categories $\llbracket \mathsf{Pth}_{\boldsymbol{\mathcal{A}}^{(2)}}\rrbracket$ and $\llbracket\mathsf{PT}_{\boldsymbol{\mathcal{A}}^{(2)}}\rrbracket$ are isomorphic.
\end{restatable}
\begin{proof}
Consider the $S$-sorted mapping $\mathrm{CH}^{\llbracket 2\rrbracket}$ introduced in Definition~\ref{DDPTQDCH}, where
\[
\mathrm{CH}^{\llbracket 2\rrbracket}
\colon
\llbracket
\mathbf{Pth}_{\boldsymbol{\mathcal{A}}^{(2)}}
\rrbracket
\mor
\llbracket\mathbf{PT}_{\boldsymbol{\mathcal{A}}^{(2)}}\rrbracket
\]

\begin{restatable}{claim}{CDIsoCatCH}
\label{CDIsoCatCH} $\mathrm{CH}^{\llbracket 2\rrbracket}$ is a $2$-functor from 
$\llbracket 
\mathsf{Pth}_{\boldsymbol{\mathcal{A}}^{(2)}}
\rrbracket$ to $
\llbracket\mathsf{PT}_{\boldsymbol{\mathcal{A}}^{(2)}}\rrbracket$.
\end{restatable}

Following Definition~\ref{DnCatSS}, we need to prove that the following equalities hold

\textsf{(1)} Let $s$ be a sort in $S$ and $\llbracket  \mathfrak{P}^{(2)}\rrbracket_{s}$ be a second-order path class in $\llbracket 
\mathsf{Pth}_{\boldsymbol{\mathcal{A}}^{(2)}}
\rrbracket_{s}$, then the following equalities hold
\begin{align*}
\mathrm{CH}^{\llbracket 2\rrbracket}_{s}\left(
\mathrm{sc}^{0\llbracket 
\mathsf{Pth}_{\boldsymbol{\mathcal{A}}^{(2)}}
\rrbracket}_{s}\left(
\llbracket  \mathfrak{P}^{(2)}\rrbracket_{s}
\right)
\right)
&=
\mathrm{sc}^{0\llbracket\mathsf{PT}_{\boldsymbol{\mathcal{A}}^{(2)}}\rrbracket}_{s}\left(
\mathrm{CH}^{\llbracket 2\rrbracket}_{s}\left(
\llbracket  \mathfrak{P}^{(2)}\rrbracket_{s}
\right)
\right);
\\
\mathrm{CH}^{\llbracket 2\rrbracket}_{s}\left(
\mathrm{tg}^{0\llbracket 
\mathsf{Pth}_{\boldsymbol{\mathcal{A}}^{(2)}}
\rrbracket}_{s}\left(
\llbracket  \mathfrak{P}^{(2)}\rrbracket_{s}
\right)
\right)
&=
\mathrm{tg}^{0\llbracket\mathsf{PT}_{\boldsymbol{\mathcal{A}}^{(2)}}\rrbracket}_{s}\left(
\mathrm{CH}^{\llbracket 2\rrbracket}_{s}\left(
\llbracket  \mathfrak{P}^{(2)}\rrbracket_{s}
\right)
\right).
\end{align*}

\textsf{(2)} Let $s$ be a sort in $S$ and $\llbracket  \mathfrak{Q}^{(2)}\rrbracket_{s}$, $\llbracket  \mathfrak{P}^{(2)}\rrbracket_{s}$ be second-order path classes in $\llbracket 
\mathsf{Pth}_{\boldsymbol{\mathcal{A}}^{(2)}}
\rrbracket_{s}$ satisfying that 
\[
\mathrm{sc}^{0\llbracket 
\mathsf{Pth}_{\boldsymbol{\mathcal{A}}^{(2)}}
\rrbracket}_{s}\left(
\llbracket  \mathfrak{Q}^{(2)}\rrbracket_{s}
\right)
=
\mathrm{tg}^{0\llbracket 
\mathsf{Pth}_{\boldsymbol{\mathcal{A}}^{(2)}}
\rrbracket}_{s}\left(
\llbracket  \mathfrak{P}^{(2)}\rrbracket_{s}
\right),
\]
then the following equality holds 
\[
\mathrm{CH}^{\llbracket 2\rrbracket}_{s}\left(
\llbracket  \mathfrak{Q}^{(2)}\rrbracket_{s}
\circ^{0\llbracket 
\mathsf{Pth}_{\boldsymbol{\mathcal{A}}^{(2)}}
\rrbracket}_{s}
\llbracket  \mathfrak{P}^{(2)}\rrbracket_{s}
\right)
=
\mathrm{CH}^{\llbracket 2\rrbracket}_{s}\left(
\llbracket  \mathfrak{Q}^{(2)}\rrbracket_{s}
\right)
\circ^{0\llbracket\mathsf{PT}_{\boldsymbol{\mathcal{A}}^{(2)}}\rrbracket}_{s}
\mathrm{CH}^{\llbracket 2\rrbracket}_{s}\left(
\llbracket  \mathfrak{P}^{(2)}\rrbracket_{s}
\right).
\]

\textsf{(3)} Let $s$ be a sort in $S$ and $\llbracket  \mathfrak{P}^{(2)}\rrbracket_{s}$ be a second-order path class in $\llbracket 
\mathsf{Pth}_{\boldsymbol{\mathcal{A}}^{(2)}}
\rrbracket_{s}$, then the following equalities hold
\begin{align*}
\mathrm{CH}^{\llbracket 2\rrbracket}_{s}\left(
\mathrm{sc}^{1\llbracket 
\mathsf{Pth}_{\boldsymbol{\mathcal{A}}^{(2)}}
\rrbracket}_{s}\left(
\llbracket  \mathfrak{P}^{(2)}\rrbracket_{s}
\right)
\right)
&=
\mathrm{sc}^{1\llbracket\mathsf{PT}_{\boldsymbol{\mathcal{A}}^{(2)}}\rrbracket}_{s}\left(
\mathrm{CH}^{\llbracket 2\rrbracket}_{s}\left(
\llbracket  \mathfrak{P}^{(2)}\rrbracket_{s}
\right)
\right);
\\
\mathrm{CH}^{\llbracket 2\rrbracket}_{s}\left(
\mathrm{tg}^{1\llbracket 
\mathsf{Pth}_{\boldsymbol{\mathcal{A}}^{(2)}}
\rrbracket}_{s}\left(
\llbracket  \mathfrak{P}^{(2)}\rrbracket_{s}
\right)
\right)
&=
\mathrm{tg}^{1\llbracket\mathsf{PT}_{\boldsymbol{\mathcal{A}}^{(2)}}\rrbracket}_{s}\left(
\mathrm{CH}^{\llbracket 2\rrbracket}_{s}\left(
\llbracket  \mathfrak{P}^{(2)}\rrbracket_{s}
\right)
\right).
\end{align*}

\textsf{(4)} Let $s$ be a sort in $S$ and $\llbracket  \mathfrak{Q}^{(2)}\rrbracket_{s}$, $\llbracket  \mathfrak{P}^{(2)}\rrbracket_{s}$ be second-order path classes in $\llbracket 
\mathsf{Pth}_{\boldsymbol{\mathcal{A}}^{(2)}}
\rrbracket_{s}$ satisfying that 
\[
\mathrm{sc}^{1\llbracket 
\mathsf{Pth}_{\boldsymbol{\mathcal{A}}^{(2)}}
\rrbracket}_{s}\left(
\llbracket  \mathfrak{Q}^{(2)}\rrbracket_{s}
\right)
=
\mathrm{tg}^{1\llbracket 
\mathsf{Pth}_{\boldsymbol{\mathcal{A}}^{(2)}}
\rrbracket}_{s}\left(
\llbracket  \mathfrak{P}^{(2)}\rrbracket_{s}
\right),
\]
then the following equality holds 
\[
\mathrm{CH}^{\llbracket 2\rrbracket}_{s}\left(
\llbracket  \mathfrak{Q}^{(2)}\rrbracket_{s}
\circ^{1\llbracket 
\mathsf{Pth}_{\boldsymbol{\mathcal{A}}^{(2)}}
\rrbracket}_{s}
\llbracket  \mathfrak{P}^{(2)}\rrbracket_{s}
\right)
=
\mathrm{CH}^{\llbracket 2\rrbracket}_{s}\left(
\llbracket  \mathfrak{Q}^{(2)}\rrbracket_{s}
\right)
\circ^{1\llbracket\mathsf{PT}_{\boldsymbol{\mathcal{A}}^{(2)}}\rrbracket}_{s}
\mathrm{CH}^{\llbracket 2\rrbracket}_{s}\left(
\llbracket  \mathfrak{P}^{(2)}\rrbracket_{s}
\right).
\]

All these equalities follows from the fact that $\mathrm{CH}^{\llbracket 2\rrbracket}$ is a many-sorted partial $\Sigma^{\boldsymbol{\mathcal{A}}^{(2)}}$-homomorphism from $\llbracket \mathbf{Pth}_{\boldsymbol{\mathcal{A}}^{(2)}}\rrbracket$ to $\llbracket\mathbf{PT}_{\boldsymbol{\mathcal{A}}^{(2)}}\rrbracket$ according to Theorem~\ref{TDIso}.

This completes the proof of Claim~\ref{CDIsoCatCH}.

Now consider the $S$-sorted mapping $\mathrm{ip}^{(\llbracket 2\rrbracket,X)@}$ introduced in Definition~\ref{DDPTQIp}, where
\[
\mathrm{ip}^{(\llbracket 2\rrbracket,X)@}
\colon
\llbracket\mathbf{PT}_{\boldsymbol{\mathcal{A}}^{(2)}}\rrbracket
\mor
\llbracket
\mathbf{Pth}_{\boldsymbol{\mathcal{A}}^{(2)}}
\rrbracket
\]

\begin{restatable}{claim}{CDIsoCatIp}
\label{CDIsoCatIp} $\mathrm{ip}^{(\llbracket 2\rrbracket,X)@}$ is a $2$-functor from  $\llbracket\mathsf{PT}_{\boldsymbol{\mathcal{A}}^{(2)}}\rrbracket$
 to $\llbracket 
\mathsf{Pth}_{\boldsymbol{\mathcal{A}}^{(2)}}
\rrbracket$.
\end{restatable}

Following Definition~\ref{DnCatSS}, we need to prove that the following equalities hold

\textsf{(1)} Let $s$ be a sort in $S$ and $\llbracket  P\rrbracket_{s}$ be a second-order path term class in $\llbracket\mathsf{PT}_{\boldsymbol{\mathcal{A}}^{(2)}}(X)_{s}\rrbracket_{s}$,
then the following equalities hold
\begin{align*}
\mathrm{ip}^{(\llbracket 2\rrbracket,X)@}_{s}\left(
\mathrm{sc}^{0\llbracket\mathsf{PT}_{\boldsymbol{\mathcal{A}}^{(2)}}\rrbracket}_{s}\left(
\llbracket  P\rrbracket_{s}
\right)
\right)
&=
\mathrm{sc}^{0\llbracket \mathsf{Pth}_{\boldsymbol{\mathcal{A}}^{(2)}}\rrbracket}_{s}\left(
\mathrm{ip}^{(\llbracket 2\rrbracket,X)@}_{s}\left(
\llbracket  P\rrbracket_{s}
\right)
\right);
\\
\mathrm{ip}^{(\llbracket 2\rrbracket,X)@}_{s}\left(
\mathrm{tg}^{0\llbracket\mathsf{PT}_{\boldsymbol{\mathcal{A}}^{(2)}}\rrbracket}_{s}\left(
\llbracket  P\rrbracket_{s}
\right)
\right)
&=
\mathrm{tg}^{0\llbracket \mathsf{Pth}_{\boldsymbol{\mathcal{A}}^{(2)}}\rrbracket}_{s}\left(
\mathrm{ip}^{(\llbracket 2\rrbracket,X)@}_{s}\left(
\llbracket  P\rrbracket_{s}
\right)
\right);
\end{align*}

\textsf{(2)} Let $s$ be a sort in $S$ and $\llbracket  Q\rrbracket_{s}$, $\llbracket  P\rrbracket_{s}$ be second-order path term classes in $\llbracket\mathsf{PT}_{\boldsymbol{\mathcal{A}}^{(2)}}(X)_{s}\rrbracket_{s}$ satisfying that 
\[
\mathrm{sc}^{0\llbracket\mathsf{PT}_{\boldsymbol{\mathcal{A}}^{(2)}}\rrbracket}_{s}\left(
\llbracket  Q\rrbracket_{s}
\right)
=
\mathrm{tg}^{0\llbracket\mathsf{PT}_{\boldsymbol{\mathcal{A}}^{(2)}}\rrbracket}_{s}\left(
\llbracket  P\rrbracket_{s}
\right),
\]
then the following equality holds
\[
\mathrm{ip}^{(\llbracket 2\rrbracket,X)@}_{s}\left(
\llbracket  Q\rrbracket_{s}
\circ^{0\llbracket\mathsf{PT}_{\boldsymbol{\mathcal{A}}^{(2)}}\rrbracket}_{s}
\llbracket  Q\rrbracket_{s}
\right)
=
\mathrm{ip}^{(\llbracket 2\rrbracket,X)@}_{s}\left(
\llbracket  Q\rrbracket_{s}
\right)
\circ^{0\llbracket \mathbf{Pth}_{\boldsymbol{\mathcal{A}}^{(2)}}\rrbracket}_{s}
\mathrm{ip}^{(\llbracket 2\rrbracket,X)@}_{s}\left(
\llbracket  P\rrbracket_{s}
\right).
\]

\textsf{(3)} Let $s$ be a sort in $S$ and $\llbracket  P\rrbracket_{s}$ be a second-order path term class in $\llbracket\mathsf{PT}_{\boldsymbol{\mathcal{A}}^{(2)}}(X)_{s}\rrbracket_{s}$, then the following equalities hold
\begin{align*}
\mathrm{ip}^{(\llbracket 2\rrbracket,X)@}_{s}\left(
\mathrm{sc}^{1\llbracket\mathsf{PT}_{\boldsymbol{\mathcal{A}}^{(2)}}\rrbracket}_{s}\left(
\llbracket  P\rrbracket_{s}
\right)
\right)
&=
\mathrm{sc}^{1\llbracket \mathsf{Pth}_{\boldsymbol{\mathcal{A}}^{(2)}}\rrbracket}_{s}\left(
\mathrm{ip}^{(\llbracket 2\rrbracket,X)@}_{s}\left(
\llbracket  P\rrbracket_{s}
\right)
\right);
\\
\mathrm{ip}^{(\llbracket 2\rrbracket,X)@}_{s}\left(
\mathrm{tg}^{1\llbracket\mathsf{PT}_{\boldsymbol{\mathcal{A}}^{(2)}}\rrbracket}_{s}\left(
\llbracket  P\rrbracket_{s}
\right)
\right)
&=
\mathrm{tg}^{1\llbracket \mathsf{Pth}_{\boldsymbol{\mathcal{A}}^{(2)}}\rrbracket}_{s}\left(
\mathrm{ip}^{(\llbracket 2\rrbracket,X)@}_{s}\left(
\llbracket  P\rrbracket_{s}
\right)
\right);
\end{align*}

\textsf{(4)} Let $s$ be a sort in $S$ and $\llbracket  Q\rrbracket_{s}$, $\llbracket  P\rrbracket_{s}$ be second-order path term classes in $\llbracket\mathsf{PT}_{\boldsymbol{\mathcal{A}}^{(2)}}(X)_{s}\rrbracket_{s}$ satisfying that 
\[
\mathrm{sc}^{1\llbracket\mathsf{PT}_{\boldsymbol{\mathcal{A}}^{(2)}}\rrbracket}_{s}\left(
\llbracket  Q\rrbracket_{s}
\right)
=
\mathrm{tg}^{1\llbracket\mathsf{PT}_{\boldsymbol{\mathcal{A}}^{(2)}}\rrbracket}_{s}\left(
\llbracket  P\rrbracket_{s}
\right),
\]
then the following equality holds
\[
\mathrm{ip}^{(\llbracket 2\rrbracket,X)@}_{s}\left(
\llbracket  Q\rrbracket_{s}
\circ^{1\llbracket\mathsf{PT}_{\boldsymbol{\mathcal{A}}^{(2)}}\rrbracket}_{s}
\llbracket  Q\rrbracket_{s}
\right)
=
\mathrm{ip}^{(\llbracket 2\rrbracket,X)@}_{s}\left(
\llbracket  Q\rrbracket_{s}
\right)
\circ^{1\llbracket \mathbf{Pth}_{\boldsymbol{\mathcal{A}}^{(2)}}\rrbracket}_{s}
\mathrm{ip}^{(\llbracket 2\rrbracket,X)@}_{s}\left(
\llbracket  P\rrbracket_{s}
\right).
\]

All these equalities follows from the fact that $\mathrm{ip}^{(\llbracket 2\rrbracket,X)@}$ is a many-sorted partial $\Sigma^{\boldsymbol{\mathcal{A}}^{(2)}}$-homomorphism from  $\llbracket\mathbf{PT}_{\boldsymbol{\mathcal{A}}^{(2)}}\rrbracket$ to $\llbracket \mathbf{Pth}_{\boldsymbol{\mathcal{A}}^{(2)}}\rrbracket$ according to Theorem~\ref{TDIso}.

This completes the proof of Claim~\ref{CDIsoCatIp}.

That the following equalities between $2$-functors hold is a direct consequence of Theorem~\ref{TDIso}.
\allowdisplaybreaks
\begin{align*}
\mathrm{ip}^{(\llbracket 2 \rrbracket,X)@}
\circ
\mathrm{CH}^{\llbracket 2\rrbracket}
&=
\mathrm{id}^{\llbracket\mathsf{Pth}_{\boldsymbol{\mathcal{A}}^{(2)}}\rrbracket};
\\
\mathrm{CH}^{\llbracket 2\rrbracket}
\circ
\mathrm{ip}^{(\llbracket 2 \rrbracket,X)@}
&=
\mathrm{id}^{\llbracket\mathsf{PT}_{\boldsymbol{\mathcal{A}}^{(2)}}\rrbracket}.
\end{align*}

This concludes the proof.
\end{proof}

As a consequence of the above results we can conclude that the $2$-categorial $\Sigma$-algebra of second-order path classes $\llbracket\mathsf{Pth}_{\boldsymbol{\mathcal{A}}^{(2)}}\rrbracket$, introduced in Definition~\ref{PDVCat}, and the $2$-categorial $\Sigma$-algebra of second-order path term classes $\llbracket\mathsf{PT}_{\boldsymbol{\mathcal{A}}^{(2)}}\rrbracket$, introduced in Definition~\ref{DDPTQDCat}, are isomorphic.

\begin{restatable}{corollary}{CDIsoCatAlg}
\label{CDIsoCatAlg} The many-sorted $2$-categorial $\Sigma$-algebras $\llbracket \mathsf{Pth}_{\boldsymbol{\mathcal{A}}^{(2)}}\rrbracket$ and $\llbracket\mathsf{PT}_{\boldsymbol{\mathcal{A}}^{(2)}}\rrbracket$ are isomorphic.
\end{restatable}
\begin{proof}
It follows from Theorem~\ref{TDIsoCat} that the $S$-sorted underlying $2$-categories $\llbracket \mathsf{Pth}_{\boldsymbol{\mathcal{A}}^{(2)}}\rrbracket$ and $\llbracket\mathsf{PT}_{\boldsymbol{\mathcal{A}}^{(2)}}\rrbracket$ are isomorphic for the pair of $2$-functors
\begin{center}
\begin{tikzpicture}
[ACliment/.style={-{To [angle'=45, length=5.75pt, width=4pt, round]}}
, scale=0.8, 
AClimentD/.style={double equal sign distance,
-implies
}
]
\node[] (PK) at (-2.9,-3) [] {$\llbracket\mathbf{Pth}_{\boldsymbol{\mathcal{A}}^{(2)}}\rrbracket$};
\node[] (PTK) at (2.9,-3) [] {$\llbracket\mathsf{PT}_{\boldsymbol{\mathcal{A}}^{(2)}}\rrbracket$};
 
 \draw[shorten >=0.2cm,shorten <=.2cm, ACliment]
 ($(-2,-3)+.09*(0,3)$) to node  [above]  {$ \mathrm{CH}^{\llbracket 2\rrbracket}$}  ($(2,-3)+.09*(0,3)$);
 \draw[shorten >=0.2cm,shorten <=.2cm, ACliment]
 ($(2,-3)-.09*(0,3)$) to node  [below]  {$\mathrm{ip}^{(\llbracket 2 \rrbracket,X)@}$}  ($(-2,-3)-.09*(0,3)$);
\end{tikzpicture}
\end{center}

All that remains to prove is the compatibility of these functors with the $\Sigma$-algebraic structure.

\begin{restatable}{claim}{CDIsoCatAlgCH}
\label{CDIsoCatAlgCH} $\mathrm{CH}^{\llbracket 2\rrbracket}$ is a $2$-categorial $\Sigma$-homomorphism from 
$\llbracket 
\mathsf{Pth}_{\boldsymbol{\mathcal{A}}^{(2)}}
\rrbracket$ to $
\llbracket\mathsf{PT}_{\boldsymbol{\mathcal{A}}^{(2)}}\rrbracket$.
\end{restatable}

Following Definition~\ref{DnCatAlg}, we need to prove that, for every pair $(\mathbf{s},s)$ in $S^{\star}\times S$, every operation $\sigma\in\Sigma_{\mathbf{s},s}$ and every  family of second-order path classes $(\llbracket \mathfrak{P}^{(2)}_{j} \rrbracket_{s_{j}})_{j\in\bb{\mathbf{s}}}$ in $\llbracket 
\mathrm{Pth}_{\boldsymbol{\mathcal{A}}^{(2)}}\rrbracket_{\mathbf{s}}$,
the following equality holds
\[
\mathrm{CH}^{\llbracket 2 \rrbracket}_{s}\left(
\sigma^{\llbracket \mathsf{Pth}_{\boldsymbol{\mathcal{A}}^{(2)}}\rrbracket}\left(\left(
\llbracket \mathfrak{P}^{(2)}_{j} \rrbracket_{s_{j}}
\right)_{j\in\bb{\mathbf{s}}}
\right)
\right)
=
\sigma^{\llbracket\mathsf{PT}_{\boldsymbol{\mathcal{A}}^{(2)}}\rrbracket}
\left(\left(
\mathrm{CH}^{\llbracket 2 \rrbracket}_{s_{j}}\left(
\llbracket \mathfrak{P}^{(2)}_{j} \rrbracket_{s_{j}}
\right)
\right)_{j\in\bb{\mathbf{s}}}
\right).
\]

This equality follows from the fact that $\mathrm{CH}^{\llbracket 2\rrbracket}$ is a many-sorted partial $\Sigma^{\boldsymbol{\mathcal{A}}^{(2)}}$-homomorphism from $\llbracket \mathbf{Pth}_{\boldsymbol{\mathcal{A}}^{(2)}}\rrbracket$ to $\llbracket\mathbf{PT}_{\boldsymbol{\mathcal{A}}^{(2)}}\rrbracket$ according to Theorem~\ref{TDIso}.

\label{******ICI}
This completes the proof of Claim~\ref{CDIsoCatAlgCH}.

\begin{restatable}{claim}{CDIsoCatAlgIp}
\label{CDIsoCatAlgIp} $\mathrm{ip}^{(\llbracket 2\rrbracket,X)@}$ is a $2$-categorial $\Sigma$-homomorphism from 
$
\llbracket\mathsf{PT}_{\boldsymbol{\mathcal{A}}^{(2)}}\rrbracket$ to 
$\llbracket 
\mathsf{Pth}_{\boldsymbol{\mathcal{A}}^{(2)}}
\rrbracket$.
\end{restatable}

Following Definition~\ref{DnCatAlg}, we need to prove that, for every pair $(\mathbf{s},s)$ in $S^{\star}\times S$, every operation $\sigma\in\Sigma_{\mathbf{s},s}$ and every  family of second-order path term classes $(\llbracket P_{j} \rrbracket_{s_{j}})_{j\in\bb{\mathbf{s}}}$ in $\llbracket\mathsf{PT}_{\boldsymbol{\mathcal{A}}^{(2)}}\rrbracket_{\mathbf{s}}$,
the following equality holds
\[
\mathrm{ip}^{(\llbracket 2\rrbracket,X)@}_{s}\left(
\sigma^{\llbracket\mathsf{PT}_{\boldsymbol{\mathcal{A}}^{(2)}}\rrbracket}\left(\left(
\llbracket P_{j} \rrbracket_{s_{j}}
\right)_{j\in\bb{\mathbf{s}}}
\right)
\right)
=
\sigma^{\llbracket 
\mathsf{Pth}_{\boldsymbol{\mathcal{A}}^{(2)}}
\rrbracket}
\left(\left(
\mathrm{ip}^{(\llbracket 2\rrbracket,X)@}_{s_{j}}\left(
\llbracket P_{j} \rrbracket_{s_{j}}
\right)
\right)_{j\in\bb{\mathbf{s}}}
\right).
\]

This equality follows from the fact that $\mathrm{ip}^{(\llbracket 2\rrbracket,X)@}$ is a many-sorted partial $\Sigma^{\boldsymbol{\mathcal{A}}^{(2)}}$-homomorphism from  $\llbracket\mathbf{PT}_{\boldsymbol{\mathcal{A}}^{(2)}}\rrbracket$ to $\llbracket \mathbf{Pth}_{\boldsymbol{\mathcal{A}}^{(2)}}\rrbracket$  according to Theorem~\ref{TDIso}.

This completes the proof of Claim~\ref{CDIsoCatAlgIp}.

This completes the proof.
\end{proof}

\section{A second-order  order isomorphism}

In this section we prove that the coproduct of $\mathrm{CH}^{\llbracket 2\rrbracket}$ and the coproduct of $\mathrm{ip}^{(\llbracket 2 \rrbracket,X)@}$ are order-preserving mappings. This will ultimately lead to prove that the Artinian partially preordered sets defined, on one hand, on Definition~\ref{DDVOrd}, on the coproduct of second-order path classes and, on the other hand, on Definition~\ref{DDPTQOrd}, on the coproduct of second-order path term classes  are order isomorphic.

We start by proving that the coproduct of the mapping $\mathrm{CH}^{\llbracket 2\rrbracket}$ is order-preserving.

\begin{restatable}{lemma}{LDIsoOrdCH}
\label{LDIsoOrdCH} 
The mapping $\coprod\mathrm{CH}^{\llbracket 2\rrbracket}$ from the set $\coprod\llbracket \mathrm{Pth}_{\boldsymbol{\mathcal{A}}^{(2)}}\rrbracket$ to the set  $\coprod\llbracket\mathrm{PT}_{\boldsymbol{\mathcal{A}}^{(2)}}\rrbracket$ that, for every sort $s$ in $S$ and every second-order path $\mathfrak{P}^{(2)}$ in $\mathrm{Pth}_{\boldsymbol{\mathcal{A}}^{(2)},s}$ sends $(\llbracket \mathfrak{P}^{(2)}\rrbracket_{s},s)$ in $\coprod\llbracket \mathrm{Pth}_{\boldsymbol{\mathcal{A}}^{(2)}}\rrbracket$ to $(\llbracket \mathrm{CH}^{(2)}_{s}(\mathfrak{P}^{(2)})\rrbracket_{s},s)$ in $\coprod\llbracket\mathrm{PT}_{\boldsymbol{\mathcal{A}}^{(2)}}\rrbracket$ determines an order-preserving mapping
\[
\textstyle
\coprod\mathrm{CH}^{\llbracket 2\rrbracket}\colon
\bigg(\coprod\llbracket \mathrm{Pth}_{\boldsymbol{\mathcal{A}}^{(2)}}\rrbracket
,\leq_{\llbracket \mathbf{Pth}_{\boldsymbol{\mathcal{A}}^{(2)}}\rrbracket}
\bigg)\mor
\left(\coprod\llbracket\mathrm{PT}_{\boldsymbol{\mathcal{A}}^{(2)}}\rrbracket,
\leq_{\llbracket \mathbf{PT}_{\boldsymbol{\mathcal{A}}^{(2)}}\rrbracket}\right).
\]
\end{restatable}
\begin{proof}
Let $s,t$ be sorts in $S$ and let us consider second-order paths $\mathfrak{Q}^{(2)}$ in $\mathrm{Pth}_{\boldsymbol{\mathcal{A}}^{(2)},t}$ and $\mathfrak{P}^{(2)}$ in $\mathrm{Pth}_{\boldsymbol{\mathcal{A}}^{(2)},s}$ satisfying that
$$
\left(\Bigl\llbracket\mathfrak{Q}^{(2)}\Bigr\rrbracket_{t},t\right)
\leq_{\llbracket \mathbf{Pth}_{\boldsymbol{\mathcal{A}}^{(2)}}\rrbracket}
\left(\Bigl\llbracket\mathfrak{P}^{(2)}\Bigr\rrbracket_{s},s\right).
$$
We have to prove that 
\[
\textstyle
\coprod\mathrm{CH}^{\llbracket 2\rrbracket}\left(\Bigl\llbracket\mathfrak{Q}^{(2)}\Bigr\rrbracket_{t},t\right)
\leq_{\llbracket \mathbf{PT}_{\boldsymbol{\mathcal{A}}^{(2)}}\rrbracket}
\coprod\mathrm{CH}^{\llbracket 2\rrbracket}\left(\Bigl\llbracket\mathfrak{P}^{(2)}\Bigr\rrbracket_{s},s\right).
\] 

Taking into account the definition of the mapping $\coprod\mathrm{CH}^{\llbracket 2\rrbracket}$, introduced in Definition~\ref{DDPTQDCH}, and the presentation of the partial order $\leq_{\llbracket \mathbf{PT}_{\boldsymbol{\mathcal{A}}^{(2)}}\rrbracket}$, introduced in  Definition~\ref{DDPTQOrd}, this is equivalent to prove the following inequality
$$
\left(
\Bigl\llbracket
\mathrm{ip}^{(2,X)@}_{t}\left(
\mathrm{CH}^{(2)}_{t}\left(
\mathfrak{Q}^{(2)}
\right)\right)\Bigr\rrbracket_{t},t\right)
\leq_{\mathbf{Pth}_{\boldsymbol{\mathcal{A}}^{(2)}}}
\left(
\Bigl\llbracket
\mathrm{ip}^{(2,X)@}_{s}\left(
\mathrm{CH}^{(2)}_{s}\left(
\mathfrak{P}^{(2)}
\right)\right)\Bigr\rrbracket_{s},s\right).
$$

Now, since $
(\llbracket\mathfrak{Q}^{(2)}\rrbracket_{t},t)\leq_{\llbracket \mathbf{Pth}_{\boldsymbol{\mathcal{A}}^{(2)}}\rrbracket}(\llbracket \mathfrak{P}^{(2)}\rrbracket_{s},s)$ we have, according to Definition~\ref{DDVOrd}, that
there exists a natural number $m\in\mathbb{N}-\{0\}$, a word $\mathbf{w}\in S^{\star}$ of length $\bb{\mathbf{w}}=m+1$ and a family of second-order paths $(\mathfrak{R}^{(2)}_{k})_{k\in\bb{\mathbf{w}}}$ in $\mathrm{Pth}_{\boldsymbol{\mathcal{A}}^{(2)},\mathbf{w}}$ such that $w_{0}=t$, $\llbracket \mathfrak{R}^{(2)}_{0}\rrbracket_{t}=\llbracket\mathfrak{Q}^{(2)}\rrbracket_{t}$, $w_{m}=s$, $\llbracket\mathfrak{R}^{(2)}_{m}\rrbracket_{s}=\llbracket\mathfrak{P}^{(2)}\rrbracket_{s}$ and, for every $k\in m$, $w_{k}=w_{k+1}$ and $\llbracket\mathfrak{R}^{(2)}_{k}\rrbracket_{w_{k}}=\llbracket\mathfrak{R}^{(2)}_{k+1}\rrbracket_{w_{m+1}}$ or $(\mathfrak{R}^{(2)}_{k}, w_{k})\leq_{\mathbf{Pth}_{\boldsymbol{\mathcal{A}}^{(2)}}} (\mathfrak{R}^{(2)}_{k+1}, w_{k+1})$.

Now, consider the family of second-order paths 
\[\left(\mathrm{ip}^{(2,X)@}_{w_{k}}\left(\mathrm{CH}^{(2)}_{w_{k}}\left(\mathfrak{R}^{(2)}_{k}
\right)\right)
\right)_{k\in\bb{\mathbf{w}}}\]
in $\mathrm{Pth}_{\boldsymbol{\mathcal{A}}^{(2)},\mathbf{w}}$. This family satisfies the following properties. 

Note that $w_{0}=t$. Moreover, the following chain of equalities holds
\allowdisplaybreaks
\begin{align*}
\left\llbracket\mathrm{ip}^{(2,X)@}_{t}\left(\mathrm{CH}^{(2)}_{t}\left(\mathfrak{R}^{(2)}_{0}
\right)\right)
\right\rrbracket_{t}
&=
\left\llbracket
\mathfrak{R}^{(2)}_{0}
\right\rrbracket_{t}
\tag{1}
\\&=
\left\llbracket
\mathfrak{Q}^{(2)}
\right\rrbracket_{t}
\tag{2}
\\&=
\left\llbracket\mathrm{ip}^{(2,X)@}_{t}\left(\mathrm{CH}^{(2)}_{t}\left(\mathfrak{Q}^{(2)}
\right)\right)
\right\rrbracket_{t}.
\tag{3}
\end{align*}

In the just stated chain of equalities, the first equality follows from Proposition~\ref{PDIpDCH}; the second equality follows by assumption; finally, the last equality follows from Proposition~\ref{PDIpDCH}.

Note also that $w_{m}=s$. Moreover, the following chain of equalities holds
\allowdisplaybreaks
\begin{align*}
\left\llbracket\mathrm{ip}^{(2,X)@}_{s}\left(\mathrm{CH}^{(2)}_{s}\left(\mathfrak{R}^{(2)}_{m}
\right)\right)
\right\rrbracket_{s}
&=
\left\llbracket
\mathfrak{R}^{(2)}_{m}
\right\rrbracket_{s}
\tag{1}
\\&=
\left\llbracket
\mathfrak{P}^{(2)}
\right\rrbracket_{s}
\tag{2}
\\&=
\left\llbracket\mathrm{ip}^{(2,X)@}_{s}\left(\mathrm{CH}^{(2)}_{s}\left(\mathfrak{P}^{(2)}
\right)\right)
\right\rrbracket_{s}.
\tag{3}
\end{align*}

In the just stated chain of equalities, the first equality follows from Proposition~\ref{PDIpDCH}; the second equality follows by assumption; finally, the last equality follows from Proposition~\ref{PDIpDCH}.

Now, for every $k\in m$, either $w_{k}=w_{k+1}$ and $\llbracket\mathfrak{R}^{(2)}_{k}\rrbracket_{w_{k}}=\llbracket\mathfrak{R}^{(2)}_{k+1}\rrbracket_{w_{m+1}}$, which will also imply that 
\[
\left\llbracket\mathrm{ip}^{(2,X)@}_{w_{k}}\left(\mathrm{CH}^{(2)}_{w_{k}}\left(\mathfrak{R}^{(2)}_{k}
\right)\right)
\right\rrbracket_{w_{k}}
=
\left\llbracket\mathrm{ip}^{(2,X)@}_{w_{k+1}}\left(\mathrm{CH}^{(2)}_{w_{k+1}}\left(\mathfrak{R}^{(2)}_{k+1}
\right)\right)
\right\rrbracket_{w_{k+1}}
,
\]
 or $(\mathfrak{R}^{(2)}_{k}, w_{k})\leq_{\mathbf{Pth}_{\boldsymbol{\mathcal{A}}^{(2)}}} (\mathfrak{R}^{(2)}_{k+1}, w_{k+1})$, which will imply, according to Corollary~\ref{CDIpDCHOrd}, that the following inequality holds
\[
\left(
\mathrm{ip}^{(2,X)@}_{w_{k}}\left(\mathrm{CH}^{(2)}_{w_{k}}\left(\mathfrak{R}^{(2)}_{k}
\right)\right),
w_{k}
\right)
\leq_{\mathbf{Pth}_{\boldsymbol{\mathcal{A}}^{(2)}}}
\left(
\mathrm{ip}^{(2,X)@}_{w_{k+1}}\left(\mathrm{CH}^{(2)}_{w_{k+1}}\left(\mathfrak{R}^{(2)}_{k+1}
\right)\right),
w_{k+1}
\right).
\]

All in all, we conclude according to Definition~\ref{DDVOrd} that
\[
\left(
\Bigl\llbracket
\mathrm{ip}^{(2,X)@}_{t}\left(
\mathrm{CH}^{(2)}_{t}\left(
\mathfrak{Q}^{(2)}
\right)\right)\Bigr\rrbracket_{t},t\right)
\leq_{\mathbf{Pth}_{\boldsymbol{\mathcal{A}}^{(2)}}}
\left(
\Bigl\llbracket
\mathrm{ip}^{(2,X)@}_{s}\left(
\mathrm{CH}^{(2)}_{s}\left(
\mathfrak{P}^{(2)}
\right)\right)\Bigr\rrbracket_{s},s\right).\]

This completes the proof.
\end{proof}

We next prove that the coproduct of the mapping $\mathrm{ip}^{(\llbracket 2 \rrbracket,X)@}$ is order-preserving.

\begin{restatable}{lemma}{LDIsoOrdIp}
\label{LDIsoOrdIp} 
The mapping $\coprod\mathrm{ip}^{(\llbracket 2 \rrbracket,X)@}$ from the set $\coprod\llbracket\mathrm{PT}_{\boldsymbol{\mathcal{A}}^{(2)}}\rrbracket$ to the set  $\coprod\llbracket \mathrm{Pth}_{\boldsymbol{\mathcal{A}}^{(2)}}\rrbracket$  that, for every sort $s$ in $S$ and every second-order path term $P$ in $\mathrm{PT}_{\boldsymbol{\mathcal{A}}^{(2)}}(X)$, sends $(\llbracket P\rrbracket_{s},s)$ in $\coprod\llbracket\mathrm{PT}_{\boldsymbol{\mathcal{A}}^{(2)}}\rrbracket$ to $(\llbracket \mathrm{ip}^{(2,X)@}_{s}(P)\rrbracket_{s},s)$ in $\coprod\llbracket \mathrm{Pth}_{\boldsymbol{\mathcal{A}}^{(2)}}\rrbracket$ determines an order-preserving mapping
\[
\textstyle
\coprod\mathrm{ip}^{(\llbracket 2 \rrbracket,X)@}\colon
\left(
\coprod\llbracket\mathrm{PT}_{\boldsymbol{\mathcal{A}}^{(2)}}\rrbracket,
\leq_{\llbracket\mathbf{PT}_{\boldsymbol{\mathcal{A}}^{(2)}}\rrbracket}\right)
\mor
\bigg(\coprod\llbracket \mathrm{Pth}_{\boldsymbol{\mathcal{A}}^{(2)}}\rrbracket
,\leq_{\llbracket \mathbf{Pth}_{\boldsymbol{\mathcal{A}}^{(2)}}\rrbracket}\bigg).
\]
\end{restatable}
\begin{proof}
Let $s,t$ be sorts in $S$ and let us consider second-order path terms $Q$ in $\mathrm{PT}_{\boldsymbol{\mathcal{A}}^{(2)},t}$ and $P$ in $\mathrm{PT}_{\boldsymbol{\mathcal{A}}^{(2)},s}$ satisfying that 
$$
\left(\llbracket
Q\rrbracket_{t},t\right)
\leq_{\llbracket\mathbf{PT}_{\boldsymbol{\mathcal{A}}^{(2)}}\rrbracket}
\left(\llbracket
P\rrbracket_{s},s\right).
$$
We have to prove that 
\[
\left(
\Bigl\llbracket \mathrm{ip}^{(2,X)@}_{t}(Q) 
\Bigr\rrbracket_{t},t
\right)
\leq_{\mathbf{Pth}_{\boldsymbol{\mathcal{A}}^{(2)}}}
\left(
\Bigl\llbracket \mathrm{ip}^{(2,X)@}_{s}(P) 
\Bigr\rrbracket_{s},s
\right).
\]

But according to Definition~\ref{DDPTQOrd}, this follows precisely from the fact that
$$
\left(\llbracket
Q\rrbracket_{t},t\right)
\leq_{\llbracket\mathbf{PT}_{\boldsymbol{\mathcal{A}}^{(2)}}\rrbracket}
\left(\llbracket
P\rrbracket_{s},s\right).
$$

This concludes  the proof.
\end{proof}

\begin{restatable}{theorem}{TDIsoOrd}
\label{TDIsoOrd} The Artinian partially preordered sets
\begin{enumerate}
\item[(i)] $(\coprod\llbracket \mathrm{Pth}_{\boldsymbol{\mathcal{A}}^{(2)}}\rrbracket
,\leq_{\llbracket \mathbf{Pth}_{\boldsymbol{\mathcal{A}}^{(2)}}\rrbracket})$; and
\item[(ii)] $(\coprod\llbracket\mathrm{PT}_{\boldsymbol{\mathcal{A}}^{(2)}}\rrbracket,
\leq_{\llbracket\mathbf{PT}_{\boldsymbol{\mathcal{A}}^{(2)}}\rrbracket})
$
\end{enumerate}
are order isomorphic.
\end{restatable}
\begin{proof}
By Lemma~\ref{LDIsoOrdCH} and Lemma~\ref{LDIsoOrdIp} the mappings $\coprod\mathrm{CH}^{\llbracket 2\rrbracket}$ and $\coprod\mathrm{ip}^{(\llbracket 2 \rrbracket,X)@}$ form a pair of mutually inverse order-preserving bijections.
\end{proof}

\chapter{
\texorpdfstring
{Freedom in the $\mathrm{QE}$-variety $\mathcal{V}(\boldsymbol{\mathcal{E}}^{\boldsymbol{\mathcal{A}}^{(2)}})$}
{Freedom in a second-order variety}
}\label{S2O}

In this chapter we define, for the second-order rewriting system $\boldsymbol{\mathcal{A}}^{(2)} = ((\mathbf{\Sigma}^{\boldsymbol{\mathcal{A}}},\boldsymbol{\mathcal{A}}),\mathcal{A}^{(2)})$, a specification $\boldsymbol{\mathcal{E}}^{\boldsymbol{\mathcal{A}}^{(2)}}$, whose defining equations 
$\mathcal{E}^{\boldsymbol{\mathcal{A}}^{(2)}}$ are $\mathrm{QE}$-equations, the $\mathrm{QE}$-variety of partial $\Sigma^{\boldsymbol{\mathcal{A}}^{(2)}}$-algebras $\mathcal{V}(\boldsymbol{\mathcal{E}}^{\boldsymbol{\mathcal{A}}^{(2)}})$ determined by $\boldsymbol{\mathcal{E}}^{\boldsymbol{\mathcal{A}}^{(2)}}$ and the category $\mathbf{PAlg}(\boldsymbol{\mathcal{E}}^{\boldsymbol{\mathcal{A}}^{(2)}})$ associated to $\mathcal{V}(\boldsymbol{\mathcal{E}}^{\boldsymbol{\mathcal{A}}^{(2)}})$. We show that the partial $\Sigma^{\boldsymbol{\mathcal{A}}^{(2)}}$-algebra of path classes $\llbracket \mathbf{Pth}_{\boldsymbol{\mathcal{A}}^{(2)}}\rrbracket$ satisfies the defining equations $\mathcal{E}^{\boldsymbol{\mathcal{A}}^{(2)}}$ of the variety $\mathcal{V}(\boldsymbol{\mathcal{E}}^{\boldsymbol{\mathcal{A}}})$, i.e., that $\llbracket\mathbf{Pth}_{\boldsymbol{\mathcal{A}}^{(2)}}\rrbracket\in \mathcal{V}(\boldsymbol{\mathcal{E}}^{\boldsymbol{\mathcal{A}}^{(2)}})$. Then we prove the fundamental result of this chapter: that the two partial $\Sigma^{\boldsymbol{\mathcal{A}}^{(2)}}$-algebras $\mathbf{T}_{\boldsymbol{\mathcal{E}}^{\boldsymbol{\mathcal{A}}^{(2)}}}(\mathbf{Pth}_{\boldsymbol{\mathcal{A}}^{(2)}})$, which is the free partial $\Sigma^{\boldsymbol{\mathcal{A}}^{(2)}}$-algebra in the category $\mathbf{PAlg}(\boldsymbol{\mathcal{E}}^{\boldsymbol{\mathcal{A}}^{(2)}})$, and $\llbracket \mathbf{Pth}_{\boldsymbol{\mathcal{A}}^{(2)}}\rrbracket$ are isomorphic. From here, as a consequence of the isomorphism of the previous chapter, we obtain that the partial $\Sigma^{\boldsymbol{\mathcal{A}}^{(2)}}$-algebra of second-order path term classes $\llbracket\mathbf{PT}_{\boldsymbol{\mathcal{A}}^{(2)}}\rrbracket$ is also a partial $\Sigma^{\boldsymbol{\mathcal{A}}^{(2)}}$-algebra in the variety $\mathcal{V}(\boldsymbol{\mathcal{E}}^{\boldsymbol{\mathcal{A}}^{(2)}})$ and isomorphic to $\mathbf{T}_{\boldsymbol{\mathcal{E}}^{\boldsymbol{\mathcal{A}}^{(2)}}}(\mathbf{Pth}_{\boldsymbol{\mathcal{A}}^{(2)}})$.

\section{
\texorpdfstring
{A $\mathrm{QE}$-variety of second-order many-sorted partial $\Sigma^{\boldsymbol{\mathcal{A}}^{(2)}}$-algebras}
{A second-order variety of partial algebras}
}

In this section we introduce, for the second-order rewriting system $\boldsymbol{\mathcal{A}}^{(2)} = ((\mathbf{\Sigma}^{\boldsymbol{\mathcal{A}}},\boldsymbol{\mathcal{A}}),\mathcal{A}^{(2)})$, the $\mathrm{QE}$-variety $\mathbf{PAlg}(\boldsymbol{\mathcal{E}}^{\boldsymbol{\mathcal{A}}^{(2)}})$ of the many-sorted partial algebras relative to the categorial signature $\Sigma^{\boldsymbol{\mathcal{A}}^{(2)}}$ and we  investigate their connections to the many-sorted partial $\Sigma^{\boldsymbol{\mathcal{A}}^{(2)}}$-algebras of second-order path classes $\llbracket\mathbf{Pth}_{\boldsymbol{\mathcal{A}}^{(2)}}\rrbracket$, introduced in Proposition~\ref{PDVDCatAlg}, and by extension, in virtue of Theorem~\ref{TDIso}, to the many-sorted partial $\Sigma^{\boldsymbol{\mathcal{A}}^{(2)}}$-algebras of second-order path term classes $\llbracket\mathbf{PT}_{\boldsymbol{\mathcal{A}}^{(2)}}\rrbracket$, introduced in Proposition~\ref{PDPTQDCatAlg}.

\begin{restatable}{definition}{DDVar}
\label{DDVar}
\index{variety!second-order!$\boldsymbol{\mathcal{E}}^{\boldsymbol{\mathcal{A}}^{(2)}}$}
\index{variety!second-order!$\mathbf{Palg}(\boldsymbol{\mathcal{E}}^{\boldsymbol{\mathcal{A}}^{(2)}})$}
For the many-sorted specification $\boldsymbol{\mathcal{A}}^{(2)} = ((\mathbf{\Sigma}^{\boldsymbol{\mathcal{A}}},\boldsymbol{\mathcal{A}}),\mathcal{A}^{(2)})$,
we will denote by $((\mathbf{\Sigma}^{\boldsymbol{\mathcal{A}}^{(2)}},V^{S}),\mathcal{E}^{\boldsymbol{\mathcal{A}}^{(2)}})$, which we write as $\boldsymbol{\mathcal{E}}^{\boldsymbol{\mathcal{A}}^{(2)}}$ for short, the many-sorted specification in which $\mathbf{\Sigma}^{\boldsymbol{\mathcal{A}}^{(2)}}$ is $(S,\Sigma^{\boldsymbol{\mathcal{A}}^{(2)}})$, $V^{S}$ a fixed $S$-sorted set with a countable infinity of variables in each coordinate, and $\mathcal{E}^{\boldsymbol{\mathcal{A}}^{(2)}}$ the subset of $\mathrm{QE}(\Sigma^{\boldsymbol{\mathcal{A}}^{(2)}})_{V^{S}}$ (which, by Definition~\ref{DQEEq}, is the set of all $\mathrm{DQEEq}$-equations with variables in $V^{S}$), consisting of the following equations:

For every $(\mathbf{s}, s)\in S^{\star}\times S$, every operation symbol $\sigma\in \Sigma_{\mathbf{s}, s}$, and every family of variables $(x_{j})_{j\in \bb{\mathbf{s}}}\in V^{S}_{\mathbf{s}}$, the operation $\sigma$ applied to the family $(x_{j})_{j\in \bb{\mathbf{s}}}$ is always defined. Formally, we have the following equation:
\allowdisplaybreaks
\begin{align*}\label{DDVarA0}\tag{A0}
\sigma((x_{j})_{j\in \bb{\mathbf{s}}})\overset{\mathrm{e}}{=}\sigma((x_{j})_{j\in \bb{\mathbf{s}}}).
\end{align*}

For every sort $s\in S$ and every variable $x\in V^{S}_{s}$, the $0$-source and $0$-target of $x$ is always defined. Formally, we have the following equations:
\allowdisplaybreaks
\begin{align*}\label{DDVarA1}\tag{A1}
\mathrm{sc}^{0}_{s}(x)&\overset{\mathrm{e}}{=}\mathrm{sc}^{0}_{s}(x);
&
\mathrm{tg}^{0}_{s}(x)&\overset{\mathrm{e}}{=}\mathrm{tg}^{0}_{s}(x).
\end{align*}

For every sort $s\in S$ and every variable $x\in V^{S}_{s}$, we have the following equations:
\allowdisplaybreaks
\begin{align*}\label{DDVarA2}
\mathrm{sc}^{0}_{s}(\mathrm{sc}^{0}_{s}(x))&\overset{\mathrm{e}}{=}\mathrm{sc}^{0}_{s}(x);
&
\mathrm{sc}^{0}_{s}(\mathrm{tg}^{0}_{s}(x))&\overset{\mathrm{e}}{=}\mathrm{tg}^{0}_{s}(x);\\
\mathrm{tg}^{0}_{s}(\mathrm{sc}^{0}_{s}(x))&\overset{\mathrm{e}}{=}\mathrm{sc}^{0}_{s}(x);
&
\mathrm{tg}^{0}_{s}(\mathrm{tg}^{0}_{s}(x))&\overset{\mathrm{e}}{=}\mathrm{tg}^{0}_{s}(x).\tag{A2}
\end{align*}
In other words, $\mathrm{sc}^{0}_{s}$ and $\mathrm{tg}^{0}_{s}$ are right zeros. In particular, $\mathrm{sc}^{0}_{s}$ and $\mathrm{tg}^{0}_{s}$ are idempotent.

For every sort $s\in S$ and every pair of variables $x,y\in V^{S}_{s}$,  $x\circ^{0}_{s} y$ is defined if and only if the $0$-target of $y$ is equal to the $0$-source of $x$. Formally, we have the following conditional equations:
\allowdisplaybreaks
\begin{align*}\label{DDVarA3}
x\circ^{0}_{s}y\overset{\mathrm{e}}{=}x\circ^{0}_{s}y \;\;
&\to \;\;
\mathrm{sc}^{0}_{s}(x)\overset{\mathrm{e}}{=}\mathrm{tg}^{0}_{s}(y);\\
\mathrm{sc}^{0}_{s}(x)\overset{\mathrm{e}}{=}\mathrm{tg}^{0}_{s}(y) \;\;
&\to \;\;
x\circ^{0}_{s}y\overset{\mathrm{e}}{=}x\circ^{0}_{s}y.\tag{A3}
\end{align*}

For every sort $s\in S$ and every pair of variables $x,y\in V^{S}_{s}$,  if $x\circ^{0}_{s} y$ is defined, then the $0$-source of $x\circ^{0}_{s}y$ is that of $y$ and the $0$-target of $x\circ^{0}_{s}y$ is that of $x$. Formally, we have the following conditional equations:
\allowdisplaybreaks
\begin{align*}\label{DDVarA4}
x\circ^{0}_{s}y\overset{\mathrm{e}}{=}x\circ^{0}_{s}y \;\;
&\to \;\;
\mathrm{sc}^{0}_{s}(x\circ^{0}_{s}y)\overset{\mathrm{e}}{=}\mathrm{sc}^{0}_{s}(y);\\
x\circ^{0}_{s}y\overset{\mathrm{e}}{=}x\circ^{0}_{s}y \;\;
&\to \;\;
\mathrm{tg}^{0}_{s}(x\circ^{0}_{s}y)\overset{\mathrm{e}}{=}\mathrm{tg}^{0}_{s}(x).\tag{A4}
\end{align*}

For every sort $s\in S$ and every variable $x\in V^{S}_{s}$, the compositions $x\circ^{0}_{s}\mathrm{sc}^{0}_{s}(x)$ and $\mathrm{tg}^{0}_{s}(x)\circ^{0}_{s}x$ are always defined and are equal to $x$, i.e., $\mathrm{sc}^{0}_{s}(x)$ is a right unit element for the $0$-composition with $x$ and $\mathrm{tg}^{0}_{s}(x)$ is a left unit element for the $0$-composition with $x$. Formally, the following equations are satisfied
\allowdisplaybreaks
\begin{align*}\label{DDVarA5}\tag{A5}
x\circ^{0}_{s}\mathrm{sc}^{0}_{s}(x)&\overset{\mathrm{e}}{=} x;
&
\mathrm{tg}^{0}_{s}(x)\circ^{0}_{s}x&\overset{\mathrm{e}}{=} x.
\end{align*}

For every sort $s\in S$ and every triple of variables $x, y, z\in V^{S}_{s}$, if the $0$-compositions $x\circ^{0}_{s}y$ and $y\circ^{0}_{s}z$ are defined, then the $0$-compositions $x\circ^{0}_{s}(y\circ^{0}_{s}z)$ and $(x\circ^{0}_{s}y)\circ^{0}_{s}z$ are defined and they are equal, i.e., the $0$-composition, when defined, is associative. Formally, we have the following conditional equation:
\allowdisplaybreaks
\begin{align*}\label{DDVarA6}\tag{A6}
(x\circ^{0}_{s}y\overset{\mathrm{e}}{=} x\circ^{0}_{s}y)
\wedge
(y\circ^{0}_{s}z\overset{\mathrm{e}}{=} y\circ^{0}_{s}z)\;\;
&\to \;\;
(x\circ^{0}_{s}y)\circ^{0}_{s}z\overset{\mathrm{e}}{=} x\circ^{0}_{s}(y\circ^{0}_{s}z).
\end{align*}

For every $(\mathbf{s}, s)\in S^{\star}\times S$, every operation symbol $\sigma\in \Sigma_{\mathbf{s}, s}$, and every family of variables $(x_{j})_{j\in \bb{\mathbf{s}}}\in V^{S}_{\mathbf{s}}$, the $0$-source of $\sigma((x_{j})_{j\in \bb{\mathbf{s}}})$ is equal to $\sigma$ applied to the family $((\mathrm{sc}^{0}_{s_{j}}(x_{j}))_{j\in\bb{\mathbf{s}}})$, and the $0$-target of $\sigma((x_{j})_{j\in \bb{\mathbf{s}}})$ is equal to $\sigma$ applied to the family $((\mathrm{tg}^{0}_{s_{j}}(x_{j}))_{j\in\bb{\mathbf{s}}})$. Formally, we have the following equations:
\allowdisplaybreaks
\begin{align*}\label{DDVarA7}
\mathrm{sc}^{0}_{s}(\sigma((x_{j})_{j\in \bb{\mathbf{s}}}))\overset{\mathrm{e}}{=}\sigma((\mathrm{sc}^{0}_{s_{j}}(x_{j}))_{j\in\bb{\mathbf{s}}});\\
\mathrm{tg}^{0}_{s}(\sigma((x_{j})_{j\in \bb{\mathbf{s}}}))\overset{\mathrm{e}}{=}\sigma((\mathrm{tg}^{0}_{s_{j}}(x_{j}))_{j\in\bb{\mathbf{s}}}).\tag{A7}
\end{align*}

For every $(\mathbf{s}, s)\in S^{\star}\times S$, every operation symbol $\sigma\in \Sigma_{\mathbf{s}, s}$, and every pair of families of variables $(x_{j})_{j\in \bb{\mathbf{s}}}, (y_{j})_{j\in \bb{\mathbf{s}}}\in V^{S}_{\mathbf{s}}$,  if, for every $j\in\bb{\mathbf{s}}$, the $0$-compositions $x_{j}\circ^{0}_{s_{j}}y_{j}$ are defined, then the $0$-composition $\sigma((x_{j})_{j\in \bb{\mathbf{s}}})\circ^{0}_{s}\sigma((y_{j})_{j\in \bb{\mathbf{s}}})$ is defined and it is equal to $\sigma$ applied to the family $(x_{j}\circ^{0}_{s_{j}}y_{j})_{j\in \bb{\mathbf{s}}}$. Formally, we have the following conditional equation:
\allowdisplaybreaks
\begin{multline*}\label{DDVarA8}
\textstyle
\bigwedge_{j\in\bb{\mathbf{s}}}
(x_{j}\circ^{0}_{s_{j}}y_{j}\overset{\mathrm{e}}{=}x_{j}\circ^{0}_{s_{j}}y_{j})
\;\;\to\;\;\\
\sigma((x_{j}\circ^{0}_{s_{j}}y_{j})_{j\in \bb{\mathbf{s}}})
\overset{\mathrm{e}}{=}
\sigma((x_{j})_{j\in \bb{\mathbf{s}}})\circ^{0}_{s}\sigma((y_{j})_{j\in \bb{\mathbf{s}}})
\tag{A8}
\end{multline*}

For every sort $s\in S$ and every rewrite rule $\mathfrak{p}\in\mathcal{A}_{s}$, $\mathfrak{p}$ is always  defined. Formally, we have the following equation:
\allowdisplaybreaks
\begin{align*}\label{DDVarA9}\tag{A9}
\mathfrak{p}\overset{\mathrm{e}}{=}\mathfrak{p}.
\end{align*}

For every sort $s\in S$ and every variable $x\in V^{S}_{s}$, the $1$-source and $1$-target of $x$ is always defined. Formally, we have the following equations:
\allowdisplaybreaks
\begin{align*}\label{DDVarB1}\tag{B1}
\mathrm{sc}^{1}_{s}(x)&\overset{\mathrm{e}}{=}\mathrm{sc}^{1}_{s}(x);
&
\mathrm{tg}^{1}_{s}(x)&\overset{\mathrm{e}}{=}\mathrm{tg}^{1}_{s}(x).
\end{align*}

For every sort $s\in S$ and every variable $x\in V^{S}_{s}$, we have the following equations:
\allowdisplaybreaks
\begin{align*}\label{DDVarB2}
\mathrm{sc}^{1}_{s}(\mathrm{sc}^{1}_{s}(x))&\overset{\mathrm{e}}{=}\mathrm{sc}^{1}_{s}(x);
&
\mathrm{sc}^{1}_{s}(\mathrm{tg}^{1}_{s}(x))&\overset{\mathrm{e}}{=}\mathrm{tg}^{1}_{s}(x);\\
\mathrm{tg}^{1}_{s}(\mathrm{sc}^{1}_{s}(x))&\overset{\mathrm{e}}{=}\mathrm{sc}^{1}_{s}(x);
&
\mathrm{tg}^{1}_{s}(\mathrm{tg}^{1}_{s}(x))&\overset{\mathrm{e}}{=}\mathrm{tg}^{1}_{s}(x).\tag{B2}
\end{align*}
In other words, $\mathrm{sc}^{1}_{s}$ and $\mathrm{tg}^{1}_{s}$ are right zeros. In particular, $\mathrm{sc}^{1}_{s}$ and $\mathrm{tg}^{1}_{s}$ are idempotent.

For every sort $s\in S$ and every pair of variables $x,y\in V^{S}_{s}$,  $x\circ^{1}_{s} y$ is defined if and only if the $1$-target of $y$ is equal to the $1$-source of $x$. Formally, we have the following conditional equations:
\allowdisplaybreaks
\begin{align*}\label{DDVarB3}
x\circ^{1}_{s}y\overset{\mathrm{e}}{=}x\circ^{1}_{s}y \;\;
&\to \;\;
\mathrm{sc}^{1}_{s}(x)\overset{\mathrm{e}}{=}\mathrm{tg}^{1}_{s}(y);\\
\mathrm{sc}^{1}_{s}(x)\overset{\mathrm{e}}{=}\mathrm{tg}^{1}_{s}(y) \;\;
&\to \;\;
x\circ^{1}_{s}y\overset{\mathrm{e}}{=}x\circ^{1}_{s}y.\tag{B3}
\end{align*}

For every sort $s\in S$ and every pair of variables $x,y\in V^{S}_{s}$,  if $x\circ^{1}_{s} y$ is defined, then the $1$-source of $x\circ^{1}_{s}y$ is that of $y$ and the $1$-target of $x\circ^{1}_{s}y$ is that of $x$. Formally, we have the following conditional equations:
\allowdisplaybreaks
\begin{align*}\label{DDVarB4}
x\circ^{1}_{s}y\overset{\mathrm{e}}{=}x\circ^{1}_{s}y \;\;
&\to \;\;
\mathrm{sc}^{1}_{s}(x\circ^{1}_{s}y)\overset{\mathrm{e}}{=}\mathrm{sc}^{1}_{s}(y);\\
x\circ^{1}_{s}y\overset{\mathrm{e}}{=}x\circ^{1}_{s}y \;\;
&\to \;\;
\mathrm{tg}^{1}_{s}(x\circ^{1}_{s}y)\overset{\mathrm{e}}{=}\mathrm{tg}^{1}_{s}(x).\tag{B4}
\end{align*}

For every sort $s\in S$ and every variable $x\in V^{S}_{s}$, the compositions $x\circ^{1}_{s}\mathrm{sc}^{1}_{s}(x)$ and $\mathrm{tg}^{1}_{s}(x)\circ^{1}_{s}x$ are always defined and are equal to $x$, i.e., $\mathrm{sc}^{1}_{s}(x)$ is a right unit element for the $1$-composition with $x$ and $\mathrm{tg}^{1}_{s}(x)$ is a left unit element for the $1$-composition with $x$. Formally, the following equations are satisfied
\allowdisplaybreaks
\begin{align*}\label{DDVarB5}\tag{B5}
x\circ^{1}_{s}\mathrm{sc}^{1}_{s}(x)&\overset{\mathrm{e}}{=} x;
&
\mathrm{tg}^{1}_{s}(x)\circ^{1}_{s}x&\overset{\mathrm{e}}{=} x.
\end{align*}

For every sort $s\in S$ and every triple of variables $x, y, z\in V^{S}_{s}$, if the $1$-compositions $x\circ^{1}_{s}y$ and $y\circ^{1}_{s}z$ are defined, then the $1$-compositions $x\circ^{1}_{s}(y\circ^{1}_{s}z)$ and $(x\circ^{1}_{s}y)\circ^{1}_{s}z$ are defined and they are equal, i.e., the $1$-composition, when defined, is associative. Formally, we have the following conditional equation:
\allowdisplaybreaks
\begin{align*}\label{DDVarB6}\tag{B6}
(x\circ^{1}_{s}y\overset{\mathrm{e}}{=} x\circ^{1}_{s}y)
\wedge
(y\circ^{1}_{s}z\overset{\mathrm{e}}{=} y\circ^{1}_{s}z)\;\;
&\to \;\;
(x\circ^{1}_{s}y)\circ^{1}_{s}z\overset{\mathrm{e}}{=} x\circ^{1}_{s}(y\circ^{1}_{s}z).
\end{align*}

For every $(\mathbf{s}, s)\in S^{\star}\times S$, every operation symbol $\sigma\in \Sigma_{\mathbf{s}, s}$, and every family of variables $(x_{j})_{j\in \bb{\mathbf{s}}}\in V^{S}_{\mathbf{s}}$, the $1$-source of $\sigma((x_{j})_{j\in \bb{\mathbf{s}}})$ is equal to $\sigma$ applied to the family $((\mathrm{sc}^{1}_{s_{j}}(x_{j}))_{j\in\bb{\mathbf{s}}})$, and the $1$-target of $\sigma((x_{j})_{j\in \bb{\mathbf{s}}})$ is equal to $\sigma$ applied to the family $((\mathrm{tg}^{1}_{s_{j}}(x_{j}))_{j\in\bb{\mathbf{s}}})$. Formally, we have the following equations:
\allowdisplaybreaks
\begin{align*}\label{DDVarB7}
\mathrm{sc}^{1}_{s}(\sigma((x_{j})_{j\in \bb{\mathbf{s}}}))\overset{\mathrm{e}}{=}\sigma((\mathrm{sc}^{1}_{s_{j}}(x_{j}))_{j\in\bb{\mathbf{s}}});\\
\mathrm{tg}^{1}_{s}(\sigma((x_{j})_{j\in \bb{\mathbf{s}}}))\overset{\mathrm{e}}{=}\sigma((\mathrm{tg}^{1}_{s_{j}}(x_{j}))_{j\in\bb{\mathbf{s}}}).\tag{B7}
\end{align*}

For every $(\mathbf{s}, s)\in S^{\star}\times S$, every operation symbol $\sigma\in \Sigma_{\mathbf{s}, s}$, and every pair of families of variables $(x_{j})_{j\in \bb{\mathbf{s}}}, (y_{j})_{j\in \bb{\mathbf{s}}}\in V^{S}_{\mathbf{s}}$,  if, for every $j\in\bb{\mathbf{s}}$, the $1$-compositions $x_{j}\circ^{1}_{s_{j}}y_{j}$ are defined, then the $1$-composition $\sigma((x_{j})_{j\in \bb{\mathbf{s}}})\circ^{1}_{s}\sigma((y_{j})_{j\in \bb{\mathbf{s}}})$ is defined and it is equal to $\sigma$ applied to the family $(x_{j}\circ^{1}_{s_{j}}y_{j})_{j\in \bb{\mathbf{s}}}$. Formally, we have the following conditional equation:
\allowdisplaybreaks
\begin{multline*}\label{DDVarB8}
\textstyle
\bigwedge_{j\in\bb{\mathbf{s}}}
(x_{j}\circ^{1}_{s_{j}}y_{j}\overset{\mathrm{e}}{=}x_{j}\circ^{1}_{s_{j}}y_{j})
\;\;\to\;\;\\
\sigma((x_{j}\circ^{1}_{s_{j}}y_{j})_{j\in \bb{\mathbf{s}}})
\overset{\mathrm{e}}{=}
\sigma((x_{j})_{j\in \bb{\mathbf{s}}})\circ^{1}_{s}\sigma((y_{j})_{j\in \bb{\mathbf{s}}})
\tag{B8}
\end{multline*}

For every sort $s\in S$ and every second-order rewrite rule $\mathfrak{p}^{(2)}\in\mathcal{A}^{(2)}_{s}$, $\mathfrak{p}^{(2)}$ is always  defined. Formally, we have the following equation:
\allowdisplaybreaks
\begin{align*}\label{DDVarB9}\tag{B9}
\mathfrak{p}^{(2)}\overset{\mathrm{e}}{=}\mathfrak{p}^{(2)}.
\end{align*}

For every sort $s\in S$ and every variable $x\in V^{S}_{s}$, the elements $\mathrm{sc}^{1}_{s}(\mathrm{sc}^{0}_{s}(x))$, $\mathrm{sc}^{0}_{s}(\mathrm{sc}^{1}_{s}(x))$ and $\mathrm{sc}^{0}_{s}(\mathrm{tg}^{1}_{s}(x))$ are always defined and are equal to $\mathrm{sc}^{0}_{s}(x)$.  Analogously, the elements $\mathrm{tg}^{1}_{s}(\mathrm{tg}^{0}_{s}(x))$, $\mathrm{tg}^{0}_{s}(\mathrm{tg}^{1}_{s}(x))$ and $\mathrm{tg}^{0}_{s}(\mathrm{sc}^{1}_{s}(x))$ are always defined and are equal to $\mathrm{tg}^{0}_{s}(x)$. Formally, the following equations are satisfied
\allowdisplaybreaks
\begin{align*}\label{DDVarAB1}\tag{AB1}
\mathrm{sc}^{1}_{s}(\mathrm{sc}^{0}_{s}(x))&\overset{\mathrm{e}}{=} \mathrm{sc}^{0}_{s}(x);
&
\mathrm{tg}^{1}_{s}(\mathrm{tg}^{0}_{s}(x))&\overset{\mathrm{e}}{=} \mathrm{tg}^{0}_{s}(x);
\\
\mathrm{sc}^{0}_{s}(\mathrm{sc}^{1}_{s}(x))&\overset{\mathrm{e}}{=} \mathrm{sc}^{0}_{s}(x);
&
\mathrm{tg}^{0}_{s}(\mathrm{tg}^{1}_{s}(x))&\overset{\mathrm{e}}{=} \mathrm{tg}^{0}_{s}(x);
\\
\mathrm{sc}^{0}_{s}(\mathrm{tg}^{1}_{s}(x))&\overset{\mathrm{e}}{=} \mathrm{sc}^{0}_{s}(x);
&
\mathrm{tg}^{0}_{s}(\mathrm{sc}^{1}_{s}(x))&\overset{\mathrm{e}}{=} \mathrm{tg}^{0}_{s}(x).
\end{align*}

For every sort $s\in S$ and every pair of variables $x,y\in V^{S}_{s}$,  if $x\circ^{0}_{s} y$ is defined, then the $1$-source of $x\circ^{0}_{s}y$ is the $0$-composition of the $1$-source of $x$ with the $1$-source of $y$ and $1$-target of $x\circ^{0}_{s}y$ is the $0$-composition of the $1$-target of $x$ with the $1$-target of $y$. Formally, we have the following conditional equations:
\allowdisplaybreaks
\begin{align*}\label{DDVarAB2}
x\circ^{0}_{s}y\overset{\mathrm{e}}{=}x\circ^{0}_{s}y \;\;
&\to \;\;
\mathrm{sc}^{1}_{s}(x\circ^{0}_{s}y)\overset{\mathrm{e}}{=}\mathrm{sc}^{1}_{s}(x)\circ^{0}_{s}\mathrm{sc}^{1}_{s}(y);\\
x\circ^{0}_{s}y\overset{\mathrm{e}}{=}x\circ^{0}_{s}y \;\;
&\to \;\;
\mathrm{tg}^{1}_{s}(x\circ^{0}_{s}y)\overset{\mathrm{e}}{=}\mathrm{tg}^{1}_{s}(x)\circ^{0}_{s}\mathrm{tg}^{1}_{s}(y).
\tag{AB2}
\end{align*}

For every sort $s\in S$ and every four variables $x,y,z,t\in V^{S}_{s}$,  if the $0$-compositions $x\circ^{0}_{s} y$ and $z\circ^{0}_{s}t$ and the $1$-compositions $x\circ^{1}_{s}z$ and $y\circ^{1}_{s}t$ are defined, then the $1$-composition of $x\circ^{0}_{s}y$ with $z\circ^{0}_{s}t$ is equal to the $0$-composition of $x\circ^{1}_{s}z$ with $y\circ^{1}_{s}t$. Formally, we have the following conditional equations:
\allowdisplaybreaks
\begin{multline*}\label{DDVarAB3}
\left(x\circ^{0}_{s}y\overset{\mathrm{e}}{=}x\circ^{0}_{s}y \right)
\wedge
\left(z\circ^{0}_{s}t\overset{\mathrm{e}}{=}z\circ^{0}_{s}t \right)
\wedge
\left(x\circ^{1}_{s}z\overset{\mathrm{e}}{=}x\circ^{1}_{s}z \right)
\wedge
\left(y\circ^{1}_{s}t\overset{\mathrm{e}}{=}y\circ^{1}_{s}t \right)
\\ \to \;\;
\left(x\circ^{0}_{s}y\right)\circ^{1}_{s}\left(z\circ^{0}_{s}t\right)
\overset{\mathrm{e}}{=}
\left(x\circ^{1}_{s}z\right)\circ^{0}_{s}\left(y\circ^{1}_{s}t\right)
.
\tag{AB3}
\end{multline*}

We will denote by $\mathbf{PAlg}(\boldsymbol{\mathcal{E}}^{\boldsymbol{\mathcal{A}}^{(2)}})$ the $\mathrm{DQEEq}$-variety determined by the many-sorted specification $\boldsymbol{\mathcal{E}}^{\boldsymbol{\mathcal{A}}^{(2)}}$.
\end{restatable}

\begin{remark}
A model of axioms $A1$--$A6$ or $B1$--$B6$ is an $S$-sorted category, see Definition~\ref{DnCat}.  A model of axioms $A1$--$A6$, $B1$--$B6$ and $AB1$--$AB3$ is an $S$-sorted $2$-category. A model of axioms $A0$--$A9$ is a $\Sigma^{\boldsymbol{\mathcal{A}}}$-algebra in $\mathbf{PAlg}(\boldsymbol{\mathcal{E}}^{\boldsymbol{\mathcal{A}}})$, see Definition~\ref{DVar}.
\end{remark}

We next prove that the many-sorted partial $\Sigma^{\boldsymbol{\mathcal{A}}^{(2)}}$-algebra $\llbracket \mathbf{Pth}_{\boldsymbol{\mathcal{A}}^{(2)}}\rrbracket$ is a model of $\boldsymbol{\mathcal{E}}^{\boldsymbol{\mathcal{A}}^{(2)}}$, i.e., is an algebra in the $\mathrm{DQEEq}$-variety $\mathbf{PAlg}(\boldsymbol{\mathcal{E}}^{\boldsymbol{\mathcal{A}}^{(2)}})$.

\begin{restatable}{proposition}{PDPthVar}
\label{PDPthVar}
$\llbracket \mathbf{Pth}_{\boldsymbol{\mathcal{A}}^{(2)}}\rrbracket$ is a many-sorted partial  $\Sigma^{\boldsymbol{\mathcal{A}}^{(2)}}$-algebra in $\mathbf{PAlg}(\boldsymbol{\mathcal{E}}^{\boldsymbol{\mathcal{A}}^{(2)}})$.
\end{restatable}

\begin{proof}
We prove that all axioms defining $\mathbf{PAlg}(\boldsymbol{\mathcal{E}}^{\boldsymbol{\mathcal{A}}^{(2)}})$ are valid in $\llbracket \mathbf{Pth}_{\boldsymbol{\mathcal{A}}^{(2)}}\rrbracket$. In this regard we recall from Proposition~\ref{PDVDCatAlg} that, since $\llbracket \mathbf{Pth}_{\boldsymbol{\mathcal{A}}^{(2)}}\rrbracket$ is a quotient of $\mathbf{Pth}_{\boldsymbol{\mathcal{A}}^{(2)}}$, the interpretation of the operation symbols in $\llbracket \mathbf{Pth}_{\boldsymbol{\mathcal{A}}^{(2)}}\rrbracket$ is defined in terms of the interpretation of the operation symbols in $\mathbf{Pth}_{\boldsymbol{\mathcal{A}}^{(2)}}$.

Axiom~\ref{DDVarA0}. Let $(\mathbf{s},s)$ be an element of $S^{\star}\times S$, $\sigma$ an operation symbol in $\Sigma_{\mathbf{s},s}$ and $(\llbracket \mathfrak{P}^{(2)}_{j}\rrbracket_{s_{j}})_{j\in\bb{\mathbf{s}}}$ be a family of second-order path classes in $\llbracket \mathrm{Pth}_{\boldsymbol{\mathcal{A}}^{(2)}}\rrbracket_{\mathbf{s}}$. Then, by Proposition~\ref{PDPthAlg}, we have that the interpretation of $\sigma$ as an operation symbol in the many-sorted partial $\Sigma^{\boldsymbol{\mathcal{A}}^{(2)}}$-algebra $\llbracket \mathbf{Pth}_{\boldsymbol{\mathcal{A}}^{(2)}}\rrbracket$ is always defined.

This proves that Axiom~\ref{DDVarA0} is valid in $\llbracket \mathbf{Pth}_{\boldsymbol{\mathcal{A}}^{(2)}}\rrbracket$.

Axiom~\ref{DDVarA1}. We recall from Proposition~\ref{PDPthCatAlg} that, for every sort $s\in S$, the $0$-source and $0$-target operations are totally defined in the many-sorted partial $\Sigma^{\boldsymbol{\mathcal{A}}^{(2)}}$-algebra $\mathbf{Pth}_{\boldsymbol{\mathcal{A}}^{(2)}}$. Let  $s$ be a sort in $S$ and $\llbracket \mathfrak{P}^{(2)}\rrbracket_{s}$ a second-order path class in $\llbracket \mathbf{Pth}_{\boldsymbol{\mathcal{A}}^{(2)}}\rrbracket_{s}$, then
\allowdisplaybreaks
\begin{align*}
\mathrm{sc}_{s}
^{0\llbracket \mathbf{Pth}_{\boldsymbol{\mathcal{A}}^{(2)}}\rrbracket}
\left(
\llbracket \mathfrak{P}^{(2)}
\rrbracket_{s}
\right)
&=
\Bigl\llbracket
\mathrm{sc}_{s}
^{0\mathbf{Pth}_{\boldsymbol{\mathcal{A}}^{(2)}}}
\left(
\mathfrak{P}^{(2)}
\right)
\Bigr\rrbracket_{s};
\\
\mathrm{tg}_{s}
^{0\llbracket \mathbf{Pth}_{\boldsymbol{\mathcal{A}}^{(2)}}\rrbracket}
\left(
\llbracket \mathfrak{P}^{(2)}
\rrbracket_{s}
\right)
&=
\Bigl\llbracket
\mathrm{tg}_{s}
^{0\mathbf{Pth}_{\boldsymbol{\mathcal{A}}^{(2)}}}
\left(
\mathfrak{P}^{(2)}
\right)
\Bigr\rrbracket_{s}.
\end{align*}

This proves that Axiom~\ref{DDVarA1} is valid in $\llbracket \mathbf{Pth}_{\boldsymbol{\mathcal{A}}^{(2)}}\rrbracket$.

Axiom~\ref{DDVarA2}. Let $s$ be a sort in $S$ and $\llbracket\mathfrak{P}\rrbracket_{s}$ a second-order path class in $\llbracket \mathbf{Pth}_{\boldsymbol{\mathcal{A}}^{(2)}}\rrbracket$. Then, by  Proposition~\ref{PDVVarA2}, we have 
\allowdisplaybreaks
\begin{align*}
\mathrm{sc}^{0\llbracket\mathbf{Pth}_{\boldsymbol{\mathcal{A}}^{(2)}}\rrbracket}_{s}\left(
\mathrm{sc}^{0\llbracket\mathbf{Pth}_{\boldsymbol{\mathcal{A}}^{(2)}}\rrbracket}_{s}\left(
\llbracket
\mathfrak{P}^{(2)}
\rrbracket_{s}
\right)
\right)
& 
=
\mathrm{sc}^{0\llbracket\mathbf{Pth}_{\boldsymbol{\mathcal{A}}^{(2)}}\rrbracket}_{s}\left(
\llbracket
\mathfrak{P}^{(2)}
\rrbracket_{s}
\right);
\\
\mathrm{sc}^{0\llbracket\mathbf{Pth}_{\boldsymbol{\mathcal{A}}^{(2)}}\rrbracket}_{s}\left(
\mathrm{tg}^{0\llbracket\mathbf{Pth}_{\boldsymbol{\mathcal{A}}^{(2)}}\rrbracket}_{s}\left(
\llbracket
\mathfrak{P}^{(2)}
\rrbracket_{s}
\right)
\right)
& 
=
\mathrm{tg}^{0\llbracket\mathbf{Pth}_{\boldsymbol{\mathcal{A}}^{(2)}}\rrbracket}_{s}\left(
\llbracket
\mathfrak{P}^{(2)}
\rrbracket_{s}
\right);
\\
\mathrm{tg}^{0\llbracket\mathbf{Pth}_{\boldsymbol{\mathcal{A}}^{(2)}}\rrbracket}_{s}\left(
\mathrm{sc}^{0\llbracket\mathbf{Pth}_{\boldsymbol{\mathcal{A}}^{(2)}}\rrbracket}_{s}\left(
\llbracket
\mathfrak{P}^{(2)}
\rrbracket_{s}
\right)
\right)
& 
=
\mathrm{sc}^{0\llbracket\mathbf{Pth}_{\boldsymbol{\mathcal{A}}^{(2)}}\rrbracket}_{s}\left(
\llbracket
\mathfrak{P}^{(2)}
\rrbracket_{s}
\right);
\\
\mathrm{tg}^{0\llbracket\mathbf{Pth}_{\boldsymbol{\mathcal{A}}^{(2)}}\rrbracket}_{s}\left(
\mathrm{tg}^{0\llbracket\mathbf{Pth}_{\boldsymbol{\mathcal{A}}^{(2)}}\rrbracket}_{s}\left(
\llbracket
\mathfrak{P}^{(2)}
\rrbracket_{s}
\right)
\right)
& 
=
\mathrm{tg}^{0\llbracket\mathbf{Pth}_{\boldsymbol{\mathcal{A}}^{(2)}}\rrbracket}_{s}\left(
\llbracket
\mathfrak{P}^{(2)}
\rrbracket_{s}
\right).
\end{align*}

This proves that Axiom~\ref{DDVarA2} is valid in $\llbracket \mathbf{Pth}_{\boldsymbol{\mathcal{A}}^{(2)}}\rrbracket$.

Axiom~\ref{DDVarA3}. Let $s$ be a sort in $S$ and $\llbracket \mathfrak{P}^{(2)}\rrbracket_{s}$, $\llbracket \mathfrak{Q}^{(2)}\rrbracket_{s}$  second-order path classes in $\llbracket\mathrm{Pth}_{\boldsymbol{\mathcal{A}}^{(2)}}\rrbracket_{s}$. By Proposition~\ref{PDVVarA3}, we can affirm that the following conditions are equivalent
\begin{enumerate}
\item[(i)] $\llbracket \mathfrak{Q}^{(2)}\rrbracket_{s}\circ^{0\llbracket\mathbf{Pth}_{\boldsymbol{\mathcal{A}}^{(2)}}\rrbracket}_{s} \llbracket \mathfrak{P}^{(2)}\rrbracket_{s}$ is defined;
\item[(ii)] $\mathrm{sc}^{0\llbracket\mathbf{Pth}_{\boldsymbol{\mathcal{A}}^{(2)}}\rrbracket}_{s}(
\llbracket\mathfrak{Q}^{(2)}\rrbracket_{s} ) = \mathrm{tg}^{0\llbracket\mathbf{Pth}_{\boldsymbol{\mathcal{A}}^{(2)}}\rrbracket}_{s}(
\llbracket\mathfrak{P}^{(2)}\rrbracket_{s} )$.
\end{enumerate}

This proves that Axiom~\ref{DDVarA3} is valid in $\llbracket \mathbf{Pth}_{\boldsymbol{\mathcal{A}}^{(2)}}\rrbracket$.

Axiom~\ref{DDVarA4}. Let $s$ be a sort in $S$ and $\llbracket \mathfrak{P}^{(2)} \rrbracket_{s}$, $\llbracket \mathfrak{Q}^{(2)} \rrbracket_{s}$ second-order path classes in $\llbracket \mathrm{Pth}_{\boldsymbol{\mathcal{A}}^{(2)}}\rrbracket_{s}$ such that the $0$-composition 
$\llbracket \mathfrak{Q}^{(2)} \rrbracket_{s}
\circ_{s}^{0\llbracket \mathbf{Pth}_{\boldsymbol{\mathcal{A}}^{(2)}}\rrbracket}
\llbracket \mathfrak{P}^{(2)} \rrbracket_{s}$
is defined.  Then, by Proposition~\ref{PDVVarA4}, we have 
\allowdisplaybreaks
\begin{align*}
\mathrm{sc}^{0\llbracket\mathbf{Pth}_{\boldsymbol{\mathcal{A}}^{(2)}}\rrbracket}_{s}\left(
\llbracket\mathfrak{Q}^{(2)}\rrbracket_{s} 
\circ^{0\llbracket\mathbf{Pth}_{\boldsymbol{\mathcal{A}}^{(2)}}\rrbracket}_{s}
\llbracket\mathfrak{P}^{(2)}\rrbracket_{s} 
\right)
&=
\mathrm{sc}^{0\llbracket\mathbf{Pth}_{\boldsymbol{\mathcal{A}}^{(2)}}\rrbracket}_{s}\left(
\llbracket\mathfrak{P}^{(2)}\rrbracket_{s} 
\right);
\\
\mathrm{tg}^{0\llbracket\mathbf{Pth}_{\boldsymbol{\mathcal{A}}^{(2)}}\rrbracket}_{s}\left(
\llbracket\mathfrak{Q}^{(2)}\rrbracket_{s} 
\circ^{0\llbracket\mathbf{Pth}_{\boldsymbol{\mathcal{A}}^{(2)}}\rrbracket}_{s}
\llbracket\mathfrak{P}^{(2)}\rrbracket_{s} 
\right)
&=
\mathrm{tg}^{0\llbracket\mathbf{Pth}_{\boldsymbol{\mathcal{A}}^{(2)}}\rrbracket}_{s}\left(
\llbracket\mathfrak{Q}^{(2)}\rrbracket_{s} 
\right).
\end{align*}

This proves that Axiom~\ref{DDVarA4} is valid in $\llbracket \mathbf{Pth}_{\boldsymbol{\mathcal{A}}^{(2)}}\rrbracket$.

Axiom~\ref{DDVarA5}. Let $s$ be a sort in $S$ and $\llbracket \mathfrak{P}^{(2)} \rrbracket_{s}$ a second-order path class in $\llbracket \mathrm{Pth}_{\boldsymbol{\mathcal{A}}^{(2)}}\rrbracket_{s}$. Then, according to Proposition~\ref{PDVVarA5}, we have that 
\allowdisplaybreaks
\begin{align*}
\llbracket\mathfrak{P}^{(2)}\rrbracket_{s} 
\circ^{0\llbracket\mathbf{Pth}_{\boldsymbol{\mathcal{A}}^{(2)}}\rrbracket}_{s}
\mathrm{sc}^{0\llbracket\mathbf{Pth}_{\boldsymbol{\mathcal{A}}^{(2)}}\rrbracket}_{s}\left(
\llbracket\mathfrak{P}^{(2)}\rrbracket_{s} 
\right)
&=
\llbracket\mathfrak{P}^{(2)}\rrbracket_{s};
\\
\mathrm{tg}^{0\llbracket\mathbf{Pth}_{\boldsymbol{\mathcal{A}}^{(2)}}\rrbracket}_{s}\left(
\llbracket\mathfrak{P}^{(2)}\rrbracket_{s} 
\right)
\circ^{0\llbracket\mathbf{Pth}_{\boldsymbol{\mathcal{A}}^{(2)}}\rrbracket}_{s}
\llbracket\mathfrak{P}^{(2)}\rrbracket_{s} 
&=
\llbracket\mathfrak{P}^{(2)}\rrbracket_{s};
\end{align*}

This proves that Axiom~\ref{DDVarA5} is valid in $\llbracket \mathbf{Pth}_{\boldsymbol{\mathcal{A}}^{(2)}}\rrbracket$.

Axiom~\ref{DDVarA6}. Let $s$ be a sort in $S$ and let $\llbracket \mathfrak{P}^{(2)} \rrbracket_{s}$, $\llbracket \mathfrak{Q}^{(2)} \rrbracket_{s}$ and $\llbracket \mathfrak{R}^{(2)} \rrbracket_{s}$ be second-order path classes in $\llbracket \mathrm{Pth}_{\boldsymbol{\mathcal{A}}^{(2)}}\rrbracket_{s}$ such that 
$
\llbracket \mathfrak{Q}^{(2)} \rrbracket_{s}
\circ^{0\llbracket \mathbf{Pth}_{\boldsymbol{\mathcal{A}}^{(2)}}\rrbracket}_{s}
\llbracket \mathfrak{P}^{(2)} \rrbracket_{s}
$ 
and 
$
\llbracket \mathfrak{R}^{(2)} \rrbracket_{s}
\circ^{0\llbracket \mathbf{Pth}_{\boldsymbol{\mathcal{A}}^{(2)}}\rrbracket}_{s}
\llbracket \mathfrak{Q}^{(2)} \rrbracket_{s}
$
are defined. Then, according to Proposition~\ref{PDVVarA6}, we have that 
\allowdisplaybreaks
\begin{multline*}
\llbracket
\mathfrak{R}^{(2)}
\rrbracket_{s}
\circ^{0\llbracket \mathbf{Pth}_{\boldsymbol{\mathcal{A}}^{(2)}}\rrbracket}_{s}
\left(
\llbracket
\mathfrak{Q}^{(2)}
\rrbracket_{s}
\circ^{0\llbracket \mathbf{Pth}_{\boldsymbol{\mathcal{A}}^{(2)}}\rrbracket}_{s}
\llbracket
\mathfrak{P}^{(2)}
\rrbracket_{s}
\right)
\\=
\left(
\llbracket
\mathfrak{R}^{(2)}
\rrbracket_{s}
\circ^{0\llbracket \mathbf{Pth}_{\boldsymbol{\mathcal{A}}^{(2)}}\rrbracket}_{s}
\llbracket
\mathfrak{Q}^{(2)}
\rrbracket_{s}
\right)
\circ^{0\llbracket \mathbf{Pth}_{\boldsymbol{\mathcal{A}}^{(2)}}\rrbracket}_{s}
\llbracket
\mathfrak{P}^{(2)}
\rrbracket_{s}.
\end{multline*}

This proves that Axiom~\ref{DDVarA6} is valid in $\llbracket \mathbf{Pth}_{\boldsymbol{\mathcal{A}}^{(2)}}\rrbracket$.

Axiom~\ref{DDVarA7}. Let $(\mathbf{s},s)$ be an element of $S^{\star}\times S$, $\sigma$ an operation symbol in $\Sigma_{\mathbf{s},s}$ and $(\llbracket \mathfrak{P}^{(2)}_{j}\rrbracket_{s_{j}})_{j\in\bb{\mathbf{s}}}$ be a family of second-order path classes in $\llbracket \mathrm{Pth}_{\boldsymbol{\mathcal{A}}^{(2)}}\rrbracket_{\mathbf{s}}$. 
Then, according to Proposition~\ref{PDVVarA7}, we have that 
\allowdisplaybreaks
\begin{multline*}
\sigma^{\llbracket \mathbf{Pth}_{\boldsymbol{\mathcal{A}}^{(2)}} \rrbracket}\left(
\left(
\mathrm{sc}^{0\llbracket \mathbf{Pth}_{\boldsymbol{\mathcal{A}}^{(2)}} \rrbracket}_{s_{j}}\left(
\Bigl\llbracket
\mathfrak{P}^{(2)}_{j}
\Bigr\rrbracket_{s_{j}}
\right)
\right)_{j\in\bb{\mathbf{s}}}
\right)
\\=
\mathrm{sc}^{0\llbracket \mathbf{Pth}_{\boldsymbol{\mathcal{A}}^{(2)}} \rrbracket}_{s}\left(
\sigma^{\llbracket \mathbf{Pth}_{\boldsymbol{\mathcal{A}}^{(2)}} \rrbracket}\left(
\left(
\Bigl\llbracket
\mathfrak{P}^{(2)}_{j}
\Bigr\rrbracket_{s_{j}}
\right)_{j\in\bb{\mathbf{s}}}
\right)
\right);
\end{multline*}
\allowdisplaybreaks
\begin{multline*}
\sigma^{\llbracket \mathbf{Pth}_{\boldsymbol{\mathcal{A}}^{(2)}} \rrbracket}\left(
\left(
\mathrm{tg}^{0\llbracket \mathbf{Pth}_{\boldsymbol{\mathcal{A}}^{(2)}} \rrbracket}_{s_{j}}\left(
\Bigl\llbracket
\mathfrak{P}^{(2)}_{j}
\Bigr\rrbracket_{s_{j}}
\right)
\right)_{j\in\bb{\mathbf{s}}}
\right)
\\=
\mathrm{tg}^{0\llbracket \mathbf{Pth}_{\boldsymbol{\mathcal{A}}^{(2)}} \rrbracket}_{s}\left(
\sigma^{\llbracket \mathbf{Pth}_{\boldsymbol{\mathcal{A}}^{(2)}} \rrbracket}\left(
\left(
\Bigl\llbracket
\mathfrak{P}^{(2)}_{j}
\Bigr\rrbracket_{s_{j}}
\right)_{j\in\bb{\mathbf{s}}}
\right)
\right).
\end{multline*}

This proves that Axiom~\ref{DDVarA7} is valid in $\llbracket \mathbf{Pth}_{\boldsymbol{\mathcal{A}}^{(2)}}\rrbracket$.

Axiom~\ref{DDVarA8}.  Let $(\mathbf{s},s)$ be an element of $S^{\star}\times S$, $\sigma\in\Sigma_{\mathbf{s},s}$ and let $(\llbracket \mathfrak{P}^{(2)}_{j}\rrbracket_{s_{j}})_{j\in\bb{\mathbf{s}}}$ and $(\llbracket \mathfrak{Q}^{(2)}_{j}\rrbracket_{s_{j}})_{j\in\bb{\mathbf{s}}}$ be two families of second-order path classes in $\llbracket \mathrm{Pth}_{\boldsymbol{\mathcal{A}}^{(2)}}\rrbracket_{\mathbf{s}}$ such that, for every $j\in\bb{\mathbf{s}}$, the $0$-compositions
$
\llbracket \mathfrak{Q}^{(2)}_{j}\rrbracket_{s_{j}}
\circ_{s_{j}}^{0\llbracket \mathbf{Pth}_{\boldsymbol{\mathcal{A}}^{(2)}}\rrbracket}
\llbracket \mathfrak{P}^{(2)}_{j}\rrbracket_{s_{j}}
$
are defined. 

Then, according to Proposition~\ref{PDVVarA8}, we have that
\allowdisplaybreaks
\begin{multline*}
\sigma^{\llbracket \mathbf{Pth}_{\boldsymbol{\mathcal{A}}^{(2)}} \rrbracket}\left(
\left(
\Bigl\llbracket
\mathfrak{Q}^{(2)}_{j}
\Bigr\rrbracket_{s_{j}}
\circ^{0\llbracket \mathbf{Pth}_{\boldsymbol{\mathcal{A}}^{(2)}} \rrbracket}_{s_{j}}
\Bigl\llbracket
\mathfrak{P}^{(2)}_{j}
\Bigr\rrbracket_{s_{j}}
\right)_{j\in\bb{\mathbf{s}}}
\right)
\\=
\sigma^{\llbracket \mathbf{Pth}_{\boldsymbol{\mathcal{A}}^{(2)}} \rrbracket}\left(
\left(
\Bigl\llbracket
\mathfrak{Q}^{(2)}_{j}
\Bigr\rrbracket_{s_{j}}
\right)_{j\in\bb{\mathbf{s}}}
\right)
\circ^{0\llbracket \mathbf{Pth}_{\boldsymbol{\mathcal{A}}^{(2)}} \rrbracket}_{s}
\sigma^{\llbracket \mathbf{Pth}_{\boldsymbol{\mathcal{A}}^{(2)}} \rrbracket}\left(
\left(
\Bigl\llbracket
\mathfrak{P}^{(2)}_{j}
\Bigr\rrbracket_{s_{j}}
\right)_{j\in\bb{\mathbf{s}}}
\right).
\end{multline*}

This proves that Axiom~\ref{DDVarA8} is valid in $\llbracket \mathbf{Pth}_{\boldsymbol{\mathcal{A}}^{(2)}}\rrbracket$.

Axiom~\ref{DDVarA9}. Let $s$ be a sort in $S$ and $\mathfrak{p}$ a rewrite rule in $\mathcal{A}_{s}$. Then, according to Proposition~\ref{PDPthCatAlg}, we have that the interpretation of $\mathfrak{p}$ as a constant in $\llbracket \mathbf{Pth}_{\boldsymbol{\mathcal{A}}^{(2)}}\rrbracket$ is given by $\mathrm{ech}^{(\llbracket 2\rrbracket,\mathcal{A})}_{s}(\mathfrak{p})$, that is, the $\Theta^{\llbracket 2\rrbracket}$-class of the $(2,[1])$-identity second-order path on the echelon determined by $\mathfrak{p}$.

This proves that Axiom~\ref{DDVarA9} is valid in $\llbracket \mathbf{Pth}_{\boldsymbol{\mathcal{A}}^{(2)}}\rrbracket$.

Axiom~\ref{DDVarB1}. We recall from Proposition~\ref{PDPthDCatAlg} that, for every sort $s\in S$, the $1$-source and $1$-target operations are totally defined in the many-sorted partial $\Sigma^{\boldsymbol{\mathcal{A}}^{(2)}}$-algebra $\mathbf{Pth}_{\boldsymbol{\mathcal{A}}^{(2)}}$. Let  $s$ be a sort in $S$ and $\llbracket \mathfrak{P}^{(2)}\rrbracket_{s}$ a second-order path class in $\llbracket \mathbf{Pth}_{\boldsymbol{\mathcal{A}}^{(2)}}\rrbracket_{s}$, then
\allowdisplaybreaks
\allowdisplaybreaks
\begin{align*}
\mathrm{sc}_{s}
^{1\llbracket \mathbf{Pth}_{\boldsymbol{\mathcal{A}}^{(2)}}\rrbracket}
\left(
\llbracket \mathfrak{P}^{(2)}
\rrbracket_{s}
\right)
&=
\Bigl\llbracket
\mathrm{sc}_{s}
^{1\mathbf{Pth}_{\boldsymbol{\mathcal{A}}^{(2)}}}
\left(
\mathfrak{P}^{(2)}
\right)
\Bigr\rrbracket_{s};
\\
\mathrm{tg}_{s}
^{1\llbracket \mathbf{Pth}_{\boldsymbol{\mathcal{A}}^{(2)}}\rrbracket}
\left(
\llbracket \mathfrak{P}^{(2)}
\rrbracket_{s}
\right)
&=
\Bigl\llbracket
\mathrm{tg}_{s}
^{1\mathbf{Pth}_{\boldsymbol{\mathcal{A}}^{(2)}}}
\left(
\mathfrak{P}^{(2)}
\right)
\Bigr\rrbracket_{s}.
\end{align*}

This proves that Axiom~\ref{DDVarB1} is valid in $\llbracket \mathbf{Pth}_{\boldsymbol{\mathcal{A}}^{(2)}}\rrbracket$.

Axiom~\ref{DDVarB2}. Let $s$ be a sort in $S$ and $\llbracket\mathfrak{P}\rrbracket_{s}$ a second-order path class in $\llbracket \mathbf{Pth}_{\boldsymbol{\mathcal{A}}^{(2)}}\rrbracket$. Then, by  Proposition~\ref{PDVVarB2}, we have 
\allowdisplaybreaks
\begin{align*}
\mathrm{sc}^{1\llbracket\mathbf{Pth}_{\boldsymbol{\mathcal{A}}^{(2)}}\rrbracket}_{s}\left(
\mathrm{sc}^{1\llbracket\mathbf{Pth}_{\boldsymbol{\mathcal{A}}^{(2)}}\rrbracket}_{s}\left(
\llbracket
\mathfrak{P}^{(2)}
\rrbracket_{s}
\right)
\right)
& 
=
\mathrm{sc}^{1\llbracket\mathbf{Pth}_{\boldsymbol{\mathcal{A}}^{(2)}}\rrbracket}_{s}\left(
\llbracket
\mathfrak{P}^{(2)}
\rrbracket_{s}
\right);
\\
\mathrm{sc}^{1\llbracket\mathbf{Pth}_{\boldsymbol{\mathcal{A}}^{(2)}}\rrbracket}_{s}\left(
\mathrm{tg}^{1\llbracket\mathbf{Pth}_{\boldsymbol{\mathcal{A}}^{(2)}}\rrbracket}_{s}\left(
\llbracket
\mathfrak{P}^{(2)}
\rrbracket_{s}
\right)
\right)
& 
=
\mathrm{tg}^{1\llbracket\mathbf{Pth}_{\boldsymbol{\mathcal{A}}^{(2)}}\rrbracket}_{s}\left(
\llbracket
\mathfrak{P}^{(2)}
\rrbracket_{s}
\right);
\\
\mathrm{tg}^{1\llbracket\mathbf{Pth}_{\boldsymbol{\mathcal{A}}^{(2)}}\rrbracket}_{s}\left(
\mathrm{sc}^{1\llbracket\mathbf{Pth}_{\boldsymbol{\mathcal{A}}^{(2)}}\rrbracket}_{s}\left(
\llbracket
\mathfrak{P}^{(2)}
\rrbracket_{s}
\right)
\right)
& 
=
\mathrm{sc}^{1\llbracket\mathbf{Pth}_{\boldsymbol{\mathcal{A}}^{(2)}}\rrbracket}_{s}\left(
\llbracket
\mathfrak{P}^{(2)}
\rrbracket_{s}
\right);
\\
\mathrm{tg}^{1\llbracket\mathbf{Pth}_{\boldsymbol{\mathcal{A}}^{(2)}}\rrbracket}_{s}\left(
\mathrm{tg}^{1\llbracket\mathbf{Pth}_{\boldsymbol{\mathcal{A}}^{(2)}}\rrbracket}_{s}\left(
\llbracket
\mathfrak{P}^{(2)}
\rrbracket_{s}
\right)
\right)
& 
=
\mathrm{tg}^{1\llbracket\mathbf{Pth}_{\boldsymbol{\mathcal{A}}^{(2)}}\rrbracket}_{s}\left(
\llbracket
\mathfrak{P}^{(2)}
\rrbracket_{s}
\right).
\end{align*}

This proves that Axiom~\ref{DDVarB2} is valid in $\llbracket \mathbf{Pth}_{\boldsymbol{\mathcal{A}}^{(2)}}\rrbracket$.

Axiom~\ref{DDVarB3}. Let $s$ be a sort in $S$ and $\llbracket \mathfrak{P}^{(2)}\rrbracket_{s}$, $\llbracket \mathfrak{Q}^{(2)}\rrbracket_{s}$  second-order path classes in $\llbracket\mathrm{Pth}_{\boldsymbol{\mathcal{A}}^{(2)}}\rrbracket_{s}$. By Proposition~\ref{PDVVarB3}, we can affirm that the following conditions are equivalent
\begin{enumerate}
\item[(i)] $\llbracket \mathfrak{Q}^{(2)}\rrbracket_{s}\circ^{1\llbracket\mathbf{Pth}_{\boldsymbol{\mathcal{A}}^{(2)}}\rrbracket}_{s} \llbracket \mathfrak{P}^{(2)}\rrbracket_{s}$ is defined;
\item[(ii)] $\mathrm{sc}^{1\llbracket\mathbf{Pth}_{\boldsymbol{\mathcal{A}}^{(2)}}\rrbracket}_{s}(
\llbracket\mathfrak{Q}^{(2)}\rrbracket_{s} ) = \mathrm{tg}^{1\llbracket\mathbf{Pth}_{\boldsymbol{\mathcal{A}}^{(2)}}\rrbracket}_{s}(
\llbracket\mathfrak{P}^{(2)}\rrbracket_{s} )$.
\end{enumerate}

This proves that Axiom~\ref{DDVarB3} is valid in $\llbracket \mathbf{Pth}_{\boldsymbol{\mathcal{A}}^{(2)}}\rrbracket$.

Axiom~\ref{DDVarB4}. Let $s$ be a sort in $S$ and $\llbracket \mathfrak{P}^{(2)} \rrbracket_{s}$, $\llbracket \mathfrak{Q}^{(2)} \rrbracket_{s}$ second-order path classes in $\llbracket \mathrm{Pth}_{\boldsymbol{\mathcal{A}}^{(2)}}\rrbracket_{s}$ such that the $1$-composition 
$\llbracket \mathfrak{Q}^{(2)} \rrbracket_{s}
\circ_{s}^{1\llbracket \mathbf{Pth}_{\boldsymbol{\mathcal{A}}^{(2)}}\rrbracket}
\llbracket \mathfrak{P}^{(2)} \rrbracket_{s}$
is defined.  Then, by Proposition~\ref{PDVVarB4}, we have 
\allowdisplaybreaks
\begin{align*}
\mathrm{sc}^{1\llbracket\mathbf{Pth}_{\boldsymbol{\mathcal{A}}^{(2)}}\rrbracket}_{s}\left(
\llbracket\mathfrak{Q}^{(2)}\rrbracket_{s} 
\circ^{1\llbracket\mathbf{Pth}_{\boldsymbol{\mathcal{A}}^{(2)}}\rrbracket}_{s}
\llbracket\mathfrak{P}^{(2)}\rrbracket_{s} 
\right)
&=
\mathrm{sc}^{1\llbracket\mathbf{Pth}_{\boldsymbol{\mathcal{A}}^{(2)}}\rrbracket}_{s}\left(
\llbracket\mathfrak{P}^{(2)}\rrbracket_{s} 
\right);
\\
\mathrm{tg}^{1\llbracket\mathbf{Pth}_{\boldsymbol{\mathcal{A}}^{(2)}}\rrbracket}_{s}\left(
\llbracket\mathfrak{Q}^{(2)}\rrbracket_{s} 
\circ^{1\llbracket\mathbf{Pth}_{\boldsymbol{\mathcal{A}}^{(2)}}\rrbracket}_{s}
\llbracket\mathfrak{P}^{(2)}\rrbracket_{s} 
\right)
&=
\mathrm{tg}^{1\llbracket\mathbf{Pth}_{\boldsymbol{\mathcal{A}}^{(2)}}\rrbracket}_{s}\left(
\llbracket\mathfrak{Q}^{(2)}\rrbracket_{s} 
\right).
\end{align*}

This proves that Axiom~\ref{DDVarB4} is valid in $\llbracket \mathbf{Pth}_{\boldsymbol{\mathcal{A}}^{(2)}}\rrbracket$.

Axiom~\ref{DDVarB5}. Let $s$ be a sort in $S$ and $\llbracket \mathfrak{P}^{(2)} \rrbracket_{s}$ a second-order path class in $\llbracket \mathrm{Pth}_{\boldsymbol{\mathcal{A}}^{(2)}}\rrbracket_{s}$. Then, according to Proposition~\ref{PDVVarB5}, we have that 
\allowdisplaybreaks
\begin{align*}
\llbracket\mathfrak{P}^{(2)}\rrbracket_{s} 
\circ^{1\llbracket\mathbf{Pth}_{\boldsymbol{\mathcal{A}}^{(2)}}\rrbracket}_{s}
\mathrm{sc}^{1\llbracket\mathbf{Pth}_{\boldsymbol{\mathcal{A}}^{(2)}}\rrbracket}_{s}\left(
\llbracket\mathfrak{P}^{(2)}\rrbracket_{s} 
\right)
&=
\llbracket\mathfrak{P}^{(2)}\rrbracket_{s};
\\
\mathrm{tg}^{1\llbracket\mathbf{Pth}_{\boldsymbol{\mathcal{A}}^{(2)}}\rrbracket}_{s}\left(
\llbracket\mathfrak{P}^{(2)}\rrbracket_{s} 
\right)
\circ^{1\llbracket\mathbf{Pth}_{\boldsymbol{\mathcal{A}}^{(2)}}\rrbracket}_{s}
\llbracket\mathfrak{P}^{(2)}\rrbracket_{s} 
&=
\llbracket\mathfrak{P}^{(2)}\rrbracket_{s};
\end{align*}

This proves that Axiom~\ref{DDVarB5} is valid in $\llbracket \mathbf{Pth}_{\boldsymbol{\mathcal{A}}^{(2)}}\rrbracket$.

Axiom~\ref{DDVarB6}. Let $s$ be a sort in $S$ and let $\llbracket \mathfrak{P}^{(2)} \rrbracket_{s}$, $\llbracket \mathfrak{Q}^{(2)} \rrbracket_{s}$ and $\llbracket \mathfrak{R}^{(2)} \rrbracket_{s}$ be second-order path classes in $\llbracket \mathrm{Pth}_{\boldsymbol{\mathcal{A}}^{(2)}}\rrbracket_{s}$ such that 
$
\llbracket \mathfrak{Q}^{(2)} \rrbracket_{s}
\circ^{1\llbracket \mathbf{Pth}_{\boldsymbol{\mathcal{A}}^{(2)}}\rrbracket}_{s}
\llbracket \mathfrak{P}^{(2)} \rrbracket_{s}
$ 
and 
$
\llbracket \mathfrak{R}^{(2)} \rrbracket_{s}
\circ^{1\llbracket \mathbf{Pth}_{\boldsymbol{\mathcal{A}}^{(2)}}\rrbracket}_{s}
\llbracket \mathfrak{Q}^{(2)} \rrbracket_{s}
$
are defined. Then, according to Proposition~\ref{PDVVarB6}, we have that 
\allowdisplaybreaks
\begin{multline*}
\llbracket
\mathfrak{R}^{(2)}
\rrbracket_{s}
\circ^{1\llbracket \mathbf{Pth}_{\boldsymbol{\mathcal{A}}^{(2)}}\rrbracket}_{s}
\left(
\llbracket
\mathfrak{Q}^{(2)}
\rrbracket_{s}
\circ^{1\llbracket \mathbf{Pth}_{\boldsymbol{\mathcal{A}}^{(2)}}\rrbracket}_{s}
\llbracket
\mathfrak{P}^{(2)}
\rrbracket_{s}
\right)
\\=
\left(
\llbracket
\mathfrak{R}^{(2)}
\rrbracket_{s}
\circ^{1\llbracket \mathbf{Pth}_{\boldsymbol{\mathcal{A}}^{(2)}}\rrbracket}_{s}
\llbracket
\mathfrak{Q}^{(2)}
\rrbracket_{s}
\right)
\circ^{1\llbracket \mathbf{Pth}_{\boldsymbol{\mathcal{A}}^{(2)}}\rrbracket}_{s}
\llbracket
\mathfrak{P}^{(2)}
\rrbracket_{s}.
\end{multline*}

This proves that Axiom~\ref{DDVarB6} is valid in $\llbracket \mathbf{Pth}_{\boldsymbol{\mathcal{A}}^{(2)}}\rrbracket$.

Axiom~\ref{DDVarB7}. Let $(\mathbf{s},s)$ be an element of $S^{\star}\times S$, $\sigma$ an operation symbol in $\Sigma_{\mathbf{s},s}$ and $(\llbracket \mathfrak{P}^{(2)}_{j}\rrbracket_{s_{j}})_{j\in\bb{\mathbf{s}}}$ be a family of second-order path classes in $\llbracket \mathrm{Pth}_{\boldsymbol{\mathcal{A}}^{(2)}}\rrbracket_{\mathbf{s}}$. 
Then, according to Proposition~\ref{PDVVarB7}, we have that 
\allowdisplaybreaks
\begin{multline*}
\sigma^{\llbracket \mathbf{Pth}_{\boldsymbol{\mathcal{A}}^{(2)}} \rrbracket}\left(
\left(
\mathrm{sc}^{1\llbracket \mathbf{Pth}_{\boldsymbol{\mathcal{A}}^{(2)}} \rrbracket}_{s_{j}}\left(
\Bigl\llbracket
\mathfrak{P}^{(2)}_{j}
\Bigr\rrbracket_{s_{j}}
\right)
\right)_{j\in\bb{\mathbf{s}}}
\right)
\\=
\mathrm{sc}^{1\llbracket \mathbf{Pth}_{\boldsymbol{\mathcal{A}}^{(2)}} \rrbracket}_{s}\left(
\sigma^{\llbracket \mathbf{Pth}_{\boldsymbol{\mathcal{A}}^{(2)}} \rrbracket}\left(
\left(
\Bigl\llbracket
\mathfrak{P}^{(2)}_{j}
\Bigr\rrbracket_{s_{j}}
\right)_{j\in\bb{\mathbf{s}}}
\right)
\right);
\end{multline*}
\allowdisplaybreaks
\begin{multline*}
\sigma^{\llbracket \mathbf{Pth}_{\boldsymbol{\mathcal{A}}^{(2)}} \rrbracket}\left(
\left(
\mathrm{tg}^{1\llbracket \mathbf{Pth}_{\boldsymbol{\mathcal{A}}^{(2)}} \rrbracket}_{s_{j}}\left(
\Bigl\llbracket
\mathfrak{P}^{(2)}_{j}
\Bigr\rrbracket_{s_{j}}
\right)
\right)_{j\in\bb{\mathbf{s}}}
\right)
\\=
\mathrm{tg}^{1\llbracket \mathbf{Pth}_{\boldsymbol{\mathcal{A}}^{(2)}} \rrbracket}_{s}\left(
\sigma^{\llbracket \mathbf{Pth}_{\boldsymbol{\mathcal{A}}^{(2)}} \rrbracket}\left(
\left(
\Bigl\llbracket
\mathfrak{P}^{(2)}_{j}
\Bigr\rrbracket_{s_{j}}
\right)_{j\in\bb{\mathbf{s}}}
\right)
\right).
\end{multline*}

This proves that Axiom~\ref{DDVarB7} is valid in $\llbracket \mathbf{Pth}_{\boldsymbol{\mathcal{A}}^{(2)}}\rrbracket$.

Axiom~\ref{DDVarB8}.  Let $(\mathbf{s},s)$ be an element of $S^{\star}\times S$, $\sigma\in\Sigma_{\mathbf{s},s}$ and let $(\llbracket \mathfrak{P}^{(2)}_{j}\rrbracket_{s_{j}})_{j\in\bb{\mathbf{s}}}$ and $(\llbracket \mathfrak{Q}^{(2)}_{j}\rrbracket_{s_{j}})_{j\in\bb{\mathbf{s}}}$ be two families of second-order path classes in $\llbracket \mathrm{Pth}_{\boldsymbol{\mathcal{A}}^{(2)}}\rrbracket_{\mathbf{s}}$ such that, for every $j\in\bb{\mathbf{s}}$, the $0$-compositions
$
\llbracket \mathfrak{Q}^{(2)}_{j}\rrbracket_{s_{j}}
\circ_{s_{j}}^{1\llbracket \mathbf{Pth}_{\boldsymbol{\mathcal{A}}^{(2)}}\rrbracket}
\llbracket \mathfrak{P}^{(2)}_{j}\rrbracket_{s_{j}}
$
are defined. 

Then, according to Proposition~\ref{PDVVarB8}, we have that
\allowdisplaybreaks
\begin{multline*}
\sigma^{\llbracket \mathbf{Pth}_{\boldsymbol{\mathcal{A}}^{(2)}} \rrbracket}\left(
\left(
\Bigl\llbracket
\mathfrak{Q}^{(2)}_{j}
\Bigr\rrbracket_{s_{j}}
\circ^{1\llbracket \mathbf{Pth}_{\boldsymbol{\mathcal{A}}^{(2)}} \rrbracket}_{s_{j}}
\Bigl\llbracket
\mathfrak{P}^{(2)}_{j}
\Bigr\rrbracket_{s_{j}}
\right)_{j\in\bb{\mathbf{s}}}
\right)
\\=
\sigma^{\llbracket \mathbf{Pth}_{\boldsymbol{\mathcal{A}}^{(2)}} \rrbracket}\left(
\left(
\Bigl\llbracket
\mathfrak{Q}^{(2)}_{j}
\Bigr\rrbracket_{s_{j}}
\right)_{j\in\bb{\mathbf{s}}}
\right)
\circ^{1\llbracket \mathbf{Pth}_{\boldsymbol{\mathcal{A}}^{(2)}} \rrbracket}_{s}
\sigma^{\llbracket \mathbf{Pth}_{\boldsymbol{\mathcal{A}}^{(2)}} \rrbracket}\left(
\left(
\Bigl\llbracket
\mathfrak{P}^{(2)}_{j}
\Bigr\rrbracket_{s_{j}}
\right)_{j\in\bb{\mathbf{s}}}
\right).
\end{multline*}

This proves that Axiom~\ref{DDVarB8} is valid in $\llbracket \mathbf{Pth}_{\boldsymbol{\mathcal{A}}^{(2)}}\rrbracket$.

Axiom~\ref{DDVarB9}. Let $s$ be a sort in $S$ and $\mathfrak{p}^{(2)}$ a second-order rewrite rule in $\mathcal{A}^{(2)}_{s}$. Then, according to Proposition~\ref{PDPthDCatAlg}, we have that the interpretation of $\mathfrak{p}^{(2)}$ as a constant in $\llbracket \mathbf{Pth}_{\boldsymbol{\mathcal{A}}^{(2)}}\rrbracket$ is given by $\mathrm{ech}^{(\llbracket 2\rrbracket,\mathcal{A}^{(2)})}_{s}(\mathfrak{p}^{(2)})$, that is, the  $\Theta^{\llbracket 2\rrbracket}$-class of second-order echelon determined by $\mathfrak{p}^{(2)}$.

This proves that Axiom~\ref{DDVarB9} is valid in $\llbracket \mathbf{Pth}_{\boldsymbol{\mathcal{A}}^{(2)}}\rrbracket$.

Axiom~\ref{DDVarAB1}.
Let $s$ be a sort in $S$ and let $\llbracket \mathfrak{P}^{(2)} \rrbracket_{s}$ be  second-order path class in $\llbracket \mathrm{Pth}_{\boldsymbol{\mathcal{A}}^{(2)}}\rrbracket_{s}$.

Then according to Proposition~\ref{PDVVarAB1}, we have that
\allowdisplaybreaks
\begin{align*}
\mathrm{sc}^{1\llbracket \mathbf{Pth}_{\boldsymbol{\mathcal{A}}^{(2)}}\rrbracket}_{s}\left(
\mathrm{sc}^{0\llbracket \mathbf{Pth}_{\boldsymbol{\mathcal{A}}^{(2)}}\rrbracket}_{s}\left(
\llbracket \mathfrak{P}^{(2)}\rrbracket_{s}
\right)\right)
&=
\mathrm{sc}^{0\llbracket \mathbf{Pth}_{\boldsymbol{\mathcal{A}}^{(2)}}\rrbracket}_{s}\left(
\llbracket \mathfrak{P}^{(2)}\rrbracket_{s}
\right);
\\
\mathrm{sc}^{0\llbracket \mathbf{Pth}_{\boldsymbol{\mathcal{A}}^{(2)}}\rrbracket}_{s}\left(
\mathrm{sc}^{1\llbracket \mathbf{Pth}_{\boldsymbol{\mathcal{A}}^{(2)}}\rrbracket}_{s}\left(
\llbracket \mathfrak{P}^{(2)}\rrbracket_{s}
\right)\right)
&=
\mathrm{sc}^{0\llbracket \mathbf{Pth}_{\boldsymbol{\mathcal{A}}^{(2)}}\rrbracket}_{s}\left(
\llbracket \mathfrak{P}^{(2)}\rrbracket_{s}
\right);
\\
\mathrm{sc}^{0\llbracket \mathbf{Pth}_{\boldsymbol{\mathcal{A}}^{(2)}}\rrbracket}_{s}\left(
\mathrm{tg}^{1\llbracket \mathbf{Pth}_{\boldsymbol{\mathcal{A}}^{(2)}}\rrbracket}_{s}\left(
\llbracket \mathfrak{P}^{(2)}\rrbracket_{s}
\right)\right)
&=
\mathrm{sc}^{0\llbracket \mathbf{Pth}_{\boldsymbol{\mathcal{A}}^{(2)}}\rrbracket}_{s}\left(
\llbracket \mathfrak{P}^{(2)}\rrbracket_{s}
\right);
\\
\mathrm{tg}^{1\llbracket \mathbf{Pth}_{\boldsymbol{\mathcal{A}}^{(2)}}\rrbracket}_{s}\left(
\mathrm{tg}^{0\llbracket \mathbf{Pth}_{\boldsymbol{\mathcal{A}}^{(2)}}\rrbracket}_{s}\left(
\llbracket \mathfrak{P}^{(2)}\rrbracket_{s}
\right)\right)
&=
\mathrm{tg}^{0\llbracket \mathbf{Pth}_{\boldsymbol{\mathcal{A}}^{(2)}}\rrbracket}_{s}\left(
\llbracket \mathfrak{P}^{(2)}\rrbracket_{s}
\right);
\\
\mathrm{tg}^{0\llbracket \mathbf{Pth}_{\boldsymbol{\mathcal{A}}^{(2)}}\rrbracket}_{s}\left(
\mathrm{tg}^{1\llbracket \mathbf{Pth}_{\boldsymbol{\mathcal{A}}^{(2)}}\rrbracket}_{s}\left(
\llbracket \mathfrak{P}^{(2)}\rrbracket_{s}
\right)\right)
&=
\mathrm{tg}^{0\llbracket \mathbf{Pth}_{\boldsymbol{\mathcal{A}}^{(2)}}\rrbracket}_{s}\left(
\llbracket \mathfrak{P}^{(2)}\rrbracket_{s}
\right);
\\
\mathrm{tg}^{0\llbracket \mathbf{Pth}_{\boldsymbol{\mathcal{A}}^{(2)}}\rrbracket}_{s}\left(
\mathrm{sc}^{1\llbracket \mathbf{Pth}_{\boldsymbol{\mathcal{A}}^{(2)}}\rrbracket}_{s}\left(
\llbracket \mathfrak{P}^{(2)}\rrbracket_{s}
\right)\right)
&=
\mathrm{tg}^{0\llbracket \mathbf{Pth}_{\boldsymbol{\mathcal{A}}^{(2)}}\rrbracket}_{s}\left(
\llbracket \mathfrak{P}^{(2)}\rrbracket_{s}
\right).
\end{align*}

This proves that Axiom~\ref{DDVarAB1} is valid in $\llbracket \mathbf{Pth}_{\boldsymbol{\mathcal{A}}^{(2)}}\rrbracket$.

Axiom~\ref{DDVarAB2}.  Let $\llbracket\mathfrak{P}^{(2)}\rrbracket_{s}$ and $\llbracket \mathfrak{Q}^{(2)}\rrbracket_{s}$ be two second-order path classes in $\llbracket \mathrm{Pth}_{\boldsymbol{\mathcal{A}}^{(2)}}\rrbracket_{s}$ such that the $0$-compositions
$
\llbracket \mathfrak{Q}^{(2)}\rrbracket_{}
\circ_{s}^{0\llbracket \mathbf{Pth}_{\boldsymbol{\mathcal{A}}^{(2)}}\rrbracket}
\llbracket \mathfrak{P}^{(2)}\rrbracket_{s}
$
is defined. 

Then, according to Proposition~\ref{PDVVarAB2}, we have that
\allowdisplaybreaks
\begin{multline*}
\mathrm{sc}^{1\llbracket \mathbf{Pth}_{\boldsymbol{\mathcal{A}}^{(2)}}\rrbracket}_{s}\left(
\llbracket \mathfrak{Q}^{(2)}\rrbracket_{s}
\circ^{0\llbracket \mathbf{Pth}_{\boldsymbol{\mathcal{A}}^{(2)}}\rrbracket}_{s}
\llbracket \mathfrak{P}^{(2)}\rrbracket_{s}
\right)
\\=
\mathrm{sc}^{1\llbracket \mathbf{Pth}_{\boldsymbol{\mathcal{A}}^{(2)}}\rrbracket}_{s}\left(
\llbracket \mathfrak{Q}^{(2)}\rrbracket_{s}
\right)
\circ^{0\llbracket \mathbf{Pth}_{\boldsymbol{\mathcal{A}}^{(2)}}\rrbracket}_{s}
\mathrm{sc}^{1\llbracket \mathbf{Pth}_{\boldsymbol{\mathcal{A}}^{(2)}}\rrbracket}_{s}\left(
\llbracket \mathfrak{P}^{(2)}\rrbracket_{s}
\right);
\end{multline*}
\allowdisplaybreaks
\begin{multline*}
\mathrm{tg}^{1\llbracket \mathbf{Pth}_{\boldsymbol{\mathcal{A}}^{(2)}}\rrbracket}_{s}\left(
\llbracket \mathfrak{Q}^{(2)}\rrbracket_{s}
\circ^{0\llbracket \mathbf{Pth}_{\boldsymbol{\mathcal{A}}^{(2)}}\rrbracket}_{s}
\llbracket \mathfrak{P}^{(2)}\rrbracket_{s}
\right)
\\=
\mathrm{tg}^{1\llbracket \mathbf{Pth}_{\boldsymbol{\mathcal{A}}^{(2)}}\rrbracket}_{s}\left(
\llbracket \mathfrak{Q}^{(2)}\rrbracket_{s}
\right)
\circ^{0\llbracket \mathbf{Pth}_{\boldsymbol{\mathcal{A}}^{(2)}}\rrbracket}_{s}
\mathrm{tg}^{1\llbracket \mathbf{Pth}_{\boldsymbol{\mathcal{A}}^{(2)}}\rrbracket}_{s}\left(
\llbracket \mathfrak{P}^{(2)}\rrbracket_{s}
\right).
\end{multline*}

This proves that Axiom~\ref{DDVarAB2} is valid in $\llbracket \mathbf{Pth}_{\boldsymbol{\mathcal{A}}^{(2)}}\rrbracket$.

Axiom~\ref{DDVarAB3}.  Let $\llbracket\mathfrak{P}^{(2)}\rrbracket_{s}$, $\llbracket\mathfrak{P}'^{(2)}\rrbracket_{s}$, $\llbracket\mathfrak{Q}^{(2)}\rrbracket_{s}$ and $\llbracket \mathfrak{Q}'^{(2)}\rrbracket_{s}$ be four second-order path classes in $\llbracket \mathrm{Pth}_{\boldsymbol{\mathcal{A}}^{(2)}}\rrbracket_{s}$ such that the $0$-compositions
$
\llbracket \mathfrak{Q}'^{(2)}\rrbracket_{}
\circ_{s}^{0\llbracket \mathbf{Pth}_{\boldsymbol{\mathcal{A}}^{(2)}}\rrbracket}
\llbracket \mathfrak{P}'^{(2)}\rrbracket_{s}
$ and $
\llbracket \mathfrak{Q}^{(2)}\rrbracket_{}
\circ_{s}^{0\llbracket \mathbf{Pth}_{\boldsymbol{\mathcal{A}}^{(2)}}\rrbracket}
\llbracket \mathfrak{P}^{(2)}\rrbracket_{s}
$
are defined and the $1$-compositions
$
\llbracket \mathfrak{Q}'^{(2)}\rrbracket_{}
\circ_{s}^{1\llbracket \mathbf{Pth}_{\boldsymbol{\mathcal{A}}^{(2)}}\rrbracket}
\llbracket \mathfrak{Q}^{(2)}\rrbracket_{s}
$ and $
\llbracket \mathfrak{P}'^{(2)}\rrbracket_{}
\circ_{s}^{1\llbracket \mathbf{Pth}_{\boldsymbol{\mathcal{A}}^{(2)}}\rrbracket}
\llbracket \mathfrak{P}^{(2)}\rrbracket_{s}
$
are defined.

Then, according to Proposition~\ref{PDVVarAB3}, we have that
\allowdisplaybreaks
\begin{multline*}
\left(
\llbracket \mathfrak{Q}'^{(2)}\rrbracket_{s}
\circ^{0\llbracket \mathbf{Pth}_{\boldsymbol{\mathcal{A}}^{(2)}}\rrbracket}_{s}
\llbracket \mathfrak{P}'^{(2)}\rrbracket_{s}
\right)
\circ^{1\llbracket \mathbf{Pth}_{\boldsymbol{\mathcal{A}}^{(2)}}\rrbracket}_{s}
\left(
\llbracket \mathfrak{Q}^{(2)}\rrbracket_{s}
\circ^{0\llbracket \mathbf{Pth}_{\boldsymbol{\mathcal{A}}^{(2)}}\rrbracket}_{s}
\llbracket \mathfrak{P}^{(2)}\rrbracket_{s}
\right)
\\=
\left(
\llbracket \mathfrak{Q}'^{(2)}\rrbracket_{s}
\circ^{1\llbracket \mathbf{Pth}_{\boldsymbol{\mathcal{A}}^{(2)}}\rrbracket}_{s}
\llbracket \mathfrak{Q}^{(2)}\rrbracket_{s}
\right)
\circ^{0\llbracket \mathbf{Pth}_{\boldsymbol{\mathcal{A}}^{(2)}}\rrbracket}_{s}
\left(
\llbracket \mathfrak{P}'^{(2)}\rrbracket_{s}
\circ^{1\llbracket \mathbf{Pth}_{\boldsymbol{\mathcal{A}}^{(2)}}\rrbracket}_{s}
\llbracket \mathfrak{P}^{(2)}\rrbracket_{s}
\right).
\end{multline*}

This proves that Axiom~\ref{DDVarAB3} is valid in $\llbracket \mathbf{Pth}_{\boldsymbol{\mathcal{A}}^{(2)}}\rrbracket$.

This completes the proof.
\end{proof}

\begin{remark} Among the axioms in $\boldsymbol{\mathcal{E}}^{\boldsymbol{\mathcal{A}}^{(2)}}$, Axioms~\ref{DDVarA5},~\ref{DDVarA6},~\ref{DDVarA8},~\ref{DDVarB8} and~\ref{DDVarAB3} are not valid in the many-sorted partial $\Sigma^{\boldsymbol{\mathcal{A}}^{(2)}}$-algebra $\mathbf{Pth}_{\boldsymbol{\mathcal{A}}^{(2)}}$.
\end{remark}

\begin{restatable}{corollary}{CDPTVar}
\label{CDPTVar} $\llbracket\mathbf{PT}_{\boldsymbol{\mathcal{A}}^{(2)}}\rrbracket$ is a many-sorted partial  $\Sigma^{\boldsymbol{\mathcal{A}}^{(2)}}$-algebra in $\mathbf{PAlg}(\boldsymbol{\mathcal{E}}^{\boldsymbol{\mathcal{A}}^{(2)}})$.
\end{restatable}
\begin{proof}
It follows from Proposition~\ref{PDPthVar} and Theorem~\ref{TDIso}.
\end{proof}

\section{
\texorpdfstring
{Freedom in $\mathcal{V}(\boldsymbol{\mathcal{E}}^{\boldsymbol{\mathcal{A}}^{(2)}})$}
{Freedom}
}

The aim of this section is to prove that the many-sorted partial $\Sigma^{\boldsymbol{\mathcal{A}}^{(2)}}$-algebra of path classes $\llbracket \mathbf{Pth}_{\boldsymbol{\mathcal{A}}^{(2)}}\rrbracket$ is a $\mathbf{PAlg}(\boldsymbol{\mathcal{E}}^{\boldsymbol{\mathcal{A}}^{(2)}})$-universal solution of the many-sorted partial $\Sigma^{\boldsymbol{\mathcal{A}}^{(2)}}$-algebra $\mathbf{Pth}_{\boldsymbol{\mathcal{A}}^{(2)}}$. In virtue of Theorem~\ref{TDIso}, the same condition will apply to the many-sorted partial $\Sigma^{\boldsymbol{\mathcal{A}}^{(2)}}$-algebra of path term classes $\llbracket\mathbf{PT}_{\boldsymbol{\mathcal{A}}^{(2)}}\rrbracket$.

The reader is advised to consult Theorem~\ref{TFreeAdj} for a description of the free many-sorted  partial $\Sigma^{\boldsymbol{\mathcal{A}}^{(2)}}$-algebra associated to a $\mathrm{DQEEq}$-variety of many-sorted partial $\Sigma^{\boldsymbol{\mathcal{A}}^{(2)}}$-algebras relative to a many-sorted partial $\Sigma^{\boldsymbol{\mathcal{A}}^{(2)}}$-algebra. 

\begin{restatable}{definition}{DDVarAbv}
\label{DDVarAbv}
\index{variety!second-order!$\mathbf{T}_{\boldsymbol{\mathcal{E}}^{\boldsymbol{\mathcal{A}}^{(2)}}}
(\mathbf{Pth}_{\boldsymbol{\mathcal{A}}^{(2)}})$}
Taking into account Theorem~\ref{TFreeAdj}, we introduce the following  abbreviations
\allowdisplaybreaks
\begin{align*}
\mathbf{Sch}
_{\boldsymbol{\mathcal{E}}^{\boldsymbol{\mathcal{A}}^{(2)}}}
(\mathbf{Pth}_{\boldsymbol{\mathcal{A}}^{(2)}})
&=
\textstyle
\bigcap
_{
f\in
\bigcup_{
\mathbf{B}\in
\mathbf{PAlg}
(
\boldsymbol{\mathcal{E}}^{\boldsymbol{\mathcal{A}}^{(2)}}
)
}
\mathrm{Hom}(
\mathbf{Pth}_{\boldsymbol{\mathcal{A}}^{(2)}}
,
\mathbf{B})
}
\mathbf{Sch}(f)
,
\\
\equiv^{\llbracket 2\rrbracket}
&=
\textstyle
\bigcap_{f\in\bigcup_{
\mathbf{B}\in
\mathbf{PAlg}
(
\boldsymbol{\mathcal{E}}^{\boldsymbol{\mathcal{A}}^{(2)}}
)
}\mathrm{Hom}(
\mathbf{Pth}_{\boldsymbol{\mathcal{A}}^{(2)}}
,
\mathbf{B}
)
}
\mathrm{SKer}(f)
,\,\text{and}
\\
\mathbf{T}_{\boldsymbol{\mathcal{E}}^{\boldsymbol{\mathcal{A}}^{(2)}}}
(\mathbf{Pth}_{\boldsymbol{\mathcal{A}}^{(2)}})
&=
\mathbf{Sch}
_{\boldsymbol{\mathcal{E}}^{\boldsymbol{\mathcal{A}}^{(2)}}}
(\mathbf{Pth}_{\boldsymbol{\mathcal{A}}^{(2)}})
/
{
\equiv^{\llbracket 2\rrbracket}
}
,
\end{align*}
where, we recall, $\mathbf{Sch}(f)$, the Schmidt closed $\mathbf{Pth}_{\boldsymbol{\mathcal{A}}^{(2)}}$-initial extension of $f$, and $\mathrm{SKer}(f)$, the Schmidt Kernel of $f$, where defined in Proposition \ref{PSch}. 
\end{restatable}

The aim of this section is to prove:
$$
\llbracket \mathbf{Pth}_{\boldsymbol{\mathcal{A}}^{(2)}}\rrbracket
\cong
\mathbf{T}_{\boldsymbol{\mathcal{E}}^{\boldsymbol{\mathcal{A}}^{(2)}}}
\left(\mathbf{Pth}_{\boldsymbol{\mathcal{A}}^{(2)}}\right).
$$

This proof will require several intermediate steps. We will start by introducing the following terminology.

\begin{restatable}{definition}{DDVarAp}
\label{DDVarAp}
We will denote by 
\begin{enumerate}
\item $\mathrm{pr}^{\equiv^{\llbracket 2\rrbracket}}$ the projection from $\mathrm{Sch}_{\boldsymbol{\mathcal{E}}^{\boldsymbol{\mathcal{A}}^{(2)}}}(\mathbf{Pth}_{\boldsymbol{\mathcal{A}}^{(2)}})$ to its quotient $\mathrm{T}_{\boldsymbol{\mathcal{E}}^{\boldsymbol{\mathcal{A}}^{(2)}}}(\mathbf{Pth}_{\boldsymbol{\mathcal{A}}^{(2)}})$;
\item  $\eta^{(\llbracket 2 \rrbracket,\mathbf{Pth}_{\boldsymbol{\mathcal{A}}^{(2)}})}$ the $S$-sorted mapping from $\mathrm{Pth}_{\boldsymbol{\mathcal{A}}^{(2)}}$ to $\mathrm{T}_{\boldsymbol{\mathcal{E}}^{\boldsymbol{\mathcal{A}}^{(2)}}}(\mathbf{Pth}_{\boldsymbol{\mathcal{A}}^{(2)}})$ given by the composition
$
\eta^{(\llbracket 2 \rrbracket,\mathbf{Pth}_{\boldsymbol{\mathcal{A}}^{(2)}})}=
\mathrm{pr}^{\equiv^{\llbracket 2\rrbracket}}\circ\eta^{(2,\mathbf{Pth}_{\boldsymbol{\mathcal{A}}^{(2)}})},
$
where $\eta^{(2,\mathbf{Pth}_{\boldsymbol{\mathcal{A}}^{(2)}})}$ is the insertion of generators, from $\mathrm{Pth}_{\boldsymbol{\mathcal{A}}^{(2)}}$ to $\mathrm{F}_{\Sigma^{\boldsymbol{\mathcal{A}}^{(2)}}}(\mathrm{Pth}_{\boldsymbol{\mathcal{A}}^{(2)}})$, introduced in Definition~\ref{DDFP}.
\end{enumerate}
\end{restatable}

\begin{remark}\label{RDVar} Let us recall from Proposition~\ref{PFreeComp} and Definition~\ref{DDFP} that $\eta^{(2,\mathbf{Pth}_{\boldsymbol{\mathcal{A}}^{(2)}})}$ is a $\Sigma^{\boldsymbol{\mathcal{A}}^{(2)}}$-homomorphism from $\mathbf{Pth}_{\boldsymbol{\mathcal{A}}^{(2)}}$ to $\mathbf{F}_{\Sigma^{\boldsymbol{\mathcal{A}}^{(2)}}}(\mathbf{Pth}_{\boldsymbol{\mathcal{A}}^{(2)}})$. Moreover, since every Schmidt algebra involved in the construction of $\mathbf{Sch}_{\boldsymbol{\mathcal{E}}^{\boldsymbol{\mathcal{A}}^{(2)}}}(\mathbf{Pth}_{\boldsymbol{\mathcal{A}}^{(2)}})$ is $\mathbf{Pth}_{\boldsymbol{\mathcal{A}}^{(2)}}$-generated, we have that $\eta^{(2,\mathbf{Pth}_{\boldsymbol{\mathcal{A}}^{(2)}})}$ corestricts to $\mathbf{Sch}_{\boldsymbol{\mathcal{E}}^{\boldsymbol{\mathcal{A}}^{(2)}}}(\mathbf{Pth}_{\boldsymbol{\mathcal{A}}^{(2)}})$.
$$
\eta^{(2,\mathbf{Pth}_{\boldsymbol{\mathcal{A}}^{(2)}})}
\colon
\mathbf{Pth}_{\boldsymbol{\mathcal{A}}^{(2)}}
\mor
\mathbf{Sch}_{\boldsymbol{\mathcal{E}}^{\boldsymbol{\mathcal{A}}^{(2)}}}(\mathbf{Pth}_{\boldsymbol{\mathcal{A}}^{(2)}}).
$$

Moreover, taking into account Definition~\ref{DDVarAbv}, the projection $\mathrm{pr}^{\equiv^{\llbracket 2\rrbracket}}$ is a $\Sigma^{\boldsymbol{\mathcal{A}}^{(2)}}$-homomorphism from $\mathbf{Sch}_{\boldsymbol{\mathcal{E}}^{\boldsymbol{\mathcal{A}}^{(2)}}}(\mathbf{Pth}_{\boldsymbol{\mathcal{A}}^{(2)}})$ to $\mathbf{T}_{\boldsymbol{\mathcal{E}}^{\boldsymbol{\mathcal{A}}^{(2)}}}(\mathbf{Pth}_{\boldsymbol{\mathcal{A}}^{(2)}})$
$$
\mathrm{pr}^{\equiv^{\llbracket 2\rrbracket}}
\colon
\mathbf{Sch}_{\boldsymbol{\mathcal{E}}^{\boldsymbol{\mathcal{A}}^{(2)}}}(\mathbf{Pth}_{\boldsymbol{\mathcal{A}}^{(2)}})
\mor
\mathbf{T}_{\boldsymbol{\mathcal{E}}^{\boldsymbol{\mathcal{A}}^{(2)}}}(\mathbf{Pth}_{\boldsymbol{\mathcal{A}}^{(2)}}).
$$

Taking into account these considerations and Definition~\ref{DDVarAp}, $\eta^{(\llbracket 2 \rrbracket,\mathbf{Pth}_{\boldsymbol{\mathcal{A}}^{(2)}})}$ is a  $\Sigma^{\boldsymbol{\mathcal{A}}^{(2)}}$-homomorphism from $\mathbf{Pth}_{\boldsymbol{\mathcal{A}}^{(2)}}$ to $\mathbf{T}_{\boldsymbol{\mathcal{E}}^{\boldsymbol{\mathcal{A}}^{(2)}}}(\mathbf{Pth}_{\boldsymbol{\mathcal{A}}^{(2)}})$
$$
\eta^{(\llbracket 2 \rrbracket,\mathbf{Pth}_{\boldsymbol{\mathcal{A}}^{(2)}})}
\colon
\mathbf{Pth}_{\boldsymbol{\mathcal{A}}^{(2)}}
\mor
\mathbf{T}_{\boldsymbol{\mathcal{E}}^{\boldsymbol{\mathcal{A}}^{(2)}}}(\mathbf{Pth}_{\boldsymbol{\mathcal{A}}^{(2)}}).
$$
\end{remark}

The following result is a consequence of Proposition~\ref{PDPthVar}. It states that the canonical projection from $\mathbf{Pth}_{\boldsymbol{\mathcal{A}}^{(2)}}$ to $\llbracket \mathbf{Pth}_{\boldsymbol{\mathcal{A}}^{(2)}}\rrbracket$ can be extended  to a $\Sigma^{\boldsymbol{\mathcal{A}}^{(2)}}$-epimorphism from $\mathbf{T}_{\boldsymbol{\mathcal{E}}^{\boldsymbol{\mathcal{A}}^{(2)}}}
(\mathbf{Pth}_{\boldsymbol{\mathcal{A}}^{(2)}})$ to $\llbracket \mathbf{Pth}_{\boldsymbol{\mathcal{A}}^{(2)}}\rrbracket$.

\begin{restatable}{corollary}{CDVarPr}
\label{CDVarPr}
For the $\Sigma^{\boldsymbol{\mathcal{A}}^{(2)}}$-epimorphism
$
\textstyle
\mathrm{pr}^{\llbracket \cdot \rrbracket}
\colon
\mathbf{Pth}_{\boldsymbol{\mathcal{A}}^{(2)}}
\mor
\llbracket \mathbf{Pth}_{\boldsymbol{\mathcal{A}}^{(2)}}\rrbracket
$, there exists a unique $\Sigma^{\boldsymbol{\mathcal{A}}^{(2)}}$-epimorphism
$$
\textstyle
\mathrm{pr}^{\llbracket \cdot \rrbracket\mathsf{p}}
\colon
\mathbf{T}_{\boldsymbol{\mathcal{E}}^{\boldsymbol{\mathcal{A}}^{(2)}}}
(\mathbf{Pth}_{\boldsymbol{\mathcal{A}}^{(2)}})
\mor
\llbracket \mathbf{Pth}_{\boldsymbol{\mathcal{A}}^{(2)}}\rrbracket
$$
such that $\mathrm{pr}^{\llbracket \cdot \rrbracket\mathsf{p}}\circ\eta^{(\llbracket 2\rrbracket,\mathbf{Pth}_{\boldsymbol{\mathcal{A}}^{(2)}})}=\mathrm{pr}^{\llbracket \cdot \rrbracket}$, i.e., such that
the diagram in Figure~\ref{FDVarPr} commutes.
\end{restatable}
\begin{figure}
\begin{tikzpicture}
[ACliment/.style={-{To [angle'=45, length=5.75pt, width=4pt, round]},
},scale=0.8]
\node[] (1) at (0,0) [] {$\mathbf{Pth}_{\boldsymbol{\mathcal{A}}^{(2)}}$};
\node[] (3) at (6,-3) []
{$\llbracket \mathbf{Pth}_{\boldsymbol{\mathcal{A}}^{(2)}}\rrbracket$};
\node[] (4) at (6,0) []
{$\mathbf{T}_{\boldsymbol{\mathcal{E}}^{\boldsymbol{\mathcal{A}}^{(2)}}}
(\mathbf{Pth}_{\boldsymbol{\mathcal{A}}^{(2)}})$};
\draw[ACliment, bend right=10]  (1) to node [below left]
{$\mathrm{pr}^{\llbracket \cdot \rrbracket}$} (3);
\draw[ACliment]  (4) to node [right]
{$\mathrm{pr}^{\llbracket \cdot \rrbracket\mathsf{p}}$} (3);
\draw[ACliment]  (1) to node [above]
{$\eta^{(\llbracket 2\rrbracket,\mathbf{Pth}_{\boldsymbol{\mathcal{A}}^{(2)}})}$} (4);
\end{tikzpicture}
\caption{The universal extension of $\mathrm{pr}^{\llbracket \cdot \rrbracket}$.}\label{FDVarPr}
\end{figure}

\begin{proof}
Let us note that $\mathrm{pr}^{\llbracket \cdot \rrbracket}$ is a $\Sigma^{\boldsymbol{\mathcal{A}}^{(2)}}$-epimorphism from
$\mathbf{Pth}_{\boldsymbol{\mathcal{A}}^{(2)}}$ to $\llbracket \mathbf{Pth}_{\boldsymbol{\mathcal{A}}^{(2)}}\rrbracket$ which, by Proposition~\ref{PDPthVar}, is a many-sorted partial $\Sigma^{\boldsymbol{\mathcal{A}}^{(2)}}$-algebra in the $\mathrm{DQEEq}$-variety $\mathbf{PAlg}(\boldsymbol{\mathcal{E}}^{\boldsymbol{\mathcal{A}}^{(2)}})$. Then, by the universal property of $\mathbf{T}_{\boldsymbol{\mathcal{E}}^{\boldsymbol{\mathcal{A}}^{(2)}}}
(\mathbf{Pth}_{\boldsymbol{\mathcal{A}}^{(2)}})$, Theorem~\ref{TFreeAdj}, there exists a unique $\Sigma^{\boldsymbol{\mathcal{A}}^{(2)}}$-homomorphism
$$
\mathrm{pr}^{\llbracket \cdot \rrbracket\mathsf{p}}
\colon
\mathbf{T}_{\boldsymbol{\mathcal{E}}^{\boldsymbol{\mathcal{A}}^{(2)}}}
(\mathbf{Pth}_{\boldsymbol{\mathcal{A}}^{(2)}})
\mor
\llbracket \mathbf{Pth}_{\boldsymbol{\mathcal{A}}^{(2)}}\rrbracket
$$
such that $\mathrm{pr}^{\llbracket \cdot \rrbracket\mathsf{p}}\circ\eta^{(\llbracket 2 \rrbracket,\mathbf{Pth}_{\boldsymbol{\mathcal{A}}^{(2)}})}=\mathrm{pr}^{\llbracket \cdot \rrbracket}$. Moreover, since $\mathrm{pr}^{\llbracket \cdot \rrbracket}$ is a $\Sigma^{\boldsymbol{\mathcal{A}}^{(2)}}$-epimorphism, we have that so is $\mathrm{pr}^{\llbracket \cdot \rrbracket\mathsf{p}}$.
\end{proof}

We now work towards the proof of the existence of an inverse of the $\Sigma^{\boldsymbol{\mathcal{A}}^{(2)}}$-homo\-mor\-phism $\mathrm{pr}^{\llbracket \cdot\rrbracket\mathsf{p}}$.

\begin{restatable}{proposition}{PDVarKer}
\label{PDVarKer}
The mapping 
$$
\mathrm{pr}^{
\equiv^{\llbracket 2\rrbracket}
}\circ\mathrm{ip}^{(2,X)@}\circ\mathrm{CH}^{(2)}
\colon
\mathbf{Pth}_{\boldsymbol{\mathcal{A}}^{(2)}}
\mor
\mathbf{T}_{\boldsymbol{\mathcal{E}}^{\boldsymbol{\mathcal{A}}^{(2)}}}
(
\mathbf{Pth}_{\boldsymbol{\mathcal{A}}^{(2)}}
)
$$
is a $\Sigma^{\boldsymbol{\mathcal{A}}^{(2)}}$-homomorphism satisfying that 
\[
\mathrm{Ker}(\mathrm{CH}^{(2)})\vee \Upsilon^{[1]}
\subseteq \mathrm{Ker}( \mathrm{pr}^{
\equiv^{\llbracket 2\rrbracket}
}\circ\mathrm{ip}^{(2,X)@}\circ\mathrm{CH}^{(2)}).
\]
\end{restatable}
\begin{proof}
Let us note that since every Schmidt algebra involved in the construction of $\mathbf{Sch}_{\boldsymbol{\mathcal{E}}^{\boldsymbol{\mathcal{A}}^{(2)}}}(\mathbf{Pth}_{\boldsymbol{\mathcal{A}}^{(2)}})$ is $\mathbf{Pth}_{\boldsymbol{\mathcal{A}}^{(2)}}$-generated and taking into account Proposition~\ref{PDIpDCH} we can correstrict the composition $\mathrm{ip}^{(2,X)@}\circ \mathrm{CH}^{(2)}$ to $\mathrm{Sch}_{\boldsymbol{\mathcal{E}}^{\boldsymbol{\mathcal{A}}^{(2)}}}(\mathbf{Pth}_{\boldsymbol{\mathcal{A}}^{(2)}})$, that is,
\[
\mathrm{ip}^{(2,X)@}\circ\mathrm{CH}^{(2)}
\colon
\mathrm{Pth}_{\boldsymbol{\mathcal{A}}^{(2)}}
\mor
\mathrm{Sch}_{\boldsymbol{\mathcal{E}}^{\boldsymbol{\mathcal{A}}^{(2)}}}(\mathbf{Pth}_{\boldsymbol{\mathcal{A}}^{(2)}}).
\]

This correstriction justifies  that  $\mathrm{pr}^{
\equiv^{\llbracket 2\rrbracket}
}\circ\mathrm{ip}^{(2,X)@}\circ\mathrm{CH}^{(2)}$ is well-defined.

We prove that $\mathrm{pr}^{
\equiv^{\llbracket 2\rrbracket}
}\circ\mathrm{ip}^{(2,X)@}\circ\mathrm{CH}^{(2)}$ is a $\Sigma^{\boldsymbol{\mathcal{A}}^{(2)}}$-homomorphism by checking the compatibility with all the different operations in $\Sigma^{\boldsymbol{\mathcal{A}}^{(2)}}$.

\textsf{Compatibility with the operations from $\Sigma$.}

Let $(\mathbf{s},s)$ be a pair in $\mathbf{S}^{\star}\times S$ and let $\sigma$ be an operation symbol in $\Sigma_{\mathbf{s},s}$, let $(\mathfrak{P}^{(2)}_{j})_{j\in\bb{\mathbf{s}}}$ be a family of second-order paths in $\mathrm{Pth}_{\boldsymbol{\mathcal{A}}^{(2)},\mathbf{s}}$. 

The following chain of equalities holds
\begin{flushleft}
$\mathrm{pr}^{\equiv^{\llbracket 2\rrbracket}}_{s}\left(
\mathrm{ip}^{(2,X)@}_{s}\left(
\mathrm{CH}^{(2)}_{s}\left(
\sigma^{\mathbf{Pth}_{\boldsymbol{\mathcal{A}}^{(2)}}}\left(
\left(
\mathfrak{P}^{(2)}_{j}
\right)_{j\in\bb{\mathbf{s}}}
\right)
\right)\right)\right)
$
\allowdisplaybreaks
\begin{align*}
\qquad&=
\mathrm{pr}^{\equiv^{\llbracket 2\rrbracket}}_{s}\left(
\sigma^{\mathbf{Pth}_{\boldsymbol{\mathcal{A}}^{(2)}}}\left(
\left(
\mathrm{ip}^{(2,X)@}_{s}\left(
\mathrm{CH}^{(2)}_{s_{j}}\left(
\mathfrak{P}^{(2)}_{j}
\right)
\right)\right)_{j\in\bb{\mathbf{s}}}
\right)\right)
\tag{1}
\\&=
\sigma^{\mathbf{T}_{\boldsymbol{\mathcal{E}}^{\boldsymbol{\mathcal{A}}^{(2)}}}
(\mathbf{Pth}_{\boldsymbol{\mathcal{A}}^{(2)}})}\left(
\left(
\mathrm{pr}^{\equiv^{\llbracket 2\rrbracket}}_{s_{j}}\left(
\mathrm{ip}^{(2,X)@}_{s}\left(
\mathrm{CH}^{(2)}_{s_{j}}\left(
\mathfrak{P}^{(2)}_{j}
\right)
\right)
\right)
\right)_{j\in\bb{\mathbf{s}}}\right).
\tag{2}
\end{align*}
\end{flushleft}

In the just stated chain of equalities, the first equality follows from Lemma~\ref{LDIpDCHSigma}; finally, the last equality follows from the fact that $\mathrm{pr}^{\equiv^{\llbracket 2\rrbracket}}$ is a $\Sigma^{\boldsymbol{\mathcal{A}}^{(2)}}$-homomorphism, according to Remark~\ref{RDVar}.

\textsf{Compatibility with the rewrite rules in $\mathcal{A}$.}

Let $s$ be a sort in $S$ and let $\mathfrak{p}$ be a rewrite rule in $\mathcal{A}_{s}$.

The following chain of equalities holds.
\allowdisplaybreaks
\begin{align*}
\mathrm{pr}^{\equiv^{\llbracket 2\rrbracket}}_{s}\left(
\mathrm{ip}^{(2,X)@}_{s}\left(
\mathrm{CH}^{(2)}_{s}\left(
\mathfrak{p}^{\mathbf{Pth}_{\boldsymbol{\mathcal{A}}^{(2)}}}
\right)\right)\right)&=
\mathrm{pr}^{\equiv^{\llbracket 2\rrbracket}}_{s}\left(
\mathfrak{p}^{\mathbf{Pth}_{\boldsymbol{\mathcal{A}}^{(2)}}}\right)
\tag{1}
\\&=
\mathfrak{p}^{\mathbf{T}_{\boldsymbol{\mathcal{E}}^{\boldsymbol{\mathcal{A}}^{(2)}}}
(\mathbf{Pth}_{\boldsymbol{\mathcal{A}}^{(2)}})}.
\tag{2}
\end{align*}

In the just stated chain of equalities, the first equality follows from Lemma~\ref{LDIpDCHEch}; finally, the last equality follows from the fact that $\mathrm{pr}^{\equiv^{\llbracket 2\rrbracket}}$ is a $\Sigma^{\boldsymbol{\mathcal{A}}^{(2)}}$-homomorphism, according to Remark~\ref{RDVar}.

\textsf{Compatibility with the $0$-source.}

Let $s$ be a sort in $S$ and let $\mathfrak{P}^{(2)}$ be a second-order path in $\mathrm{Pth}_{\boldsymbol{\mathcal{A}}^{(2)},s}$. 

The following chain of equalities holds
\begin{flushleft}
$\mathrm{pr}^{\equiv^{\llbracket 2\rrbracket}}_{s}\left(
\mathrm{ip}^{(2,X)@}_{s}\left(
\mathrm{CH}^{(2)}_{s}\left(
\mathrm{sc}^{0\mathbf{Pth}_{\boldsymbol{\mathcal{A}}^{(2)}}}_{s}\left(
\mathfrak{P}^{(2)}
\right)
\right)\right)\right)
$
\allowdisplaybreaks
\begin{align*}
\qquad&=
\mathrm{pr}^{\equiv^{\llbracket 2\rrbracket}}_{s}\left(
\mathrm{sc}^{0\mathbf{Pth}_{\boldsymbol{\mathcal{A}}^{(2)}}}_{s}\left(
\mathrm{ip}^{(2,X)@}_{s}\left(
\mathrm{CH}^{(2)}_{s}\left(
\mathfrak{P}^{(2)}
\right)
\right)\right)\right)
\tag{1}
\\&=
\mathrm{sc}^{0\mathbf{T}_{\boldsymbol{\mathcal{E}}^{\boldsymbol{\mathcal{A}}^{(2)}}}
(\mathbf{Pth}_{\boldsymbol{\mathcal{A}}^{(2)}})}_{s}\left(
\mathrm{pr}^{\equiv^{\llbracket 2\rrbracket}}_{s}\left(
\mathrm{ip}^{(2,X)@}_{s}\left(
\mathrm{CH}^{(2)}_{s}\left(
\mathfrak{P}^{(2)}
\right)
\right)\right)\right).
\tag{2}
\end{align*}
\end{flushleft}

In the just stated chain of equalities, the first equality follows from Lemma~\ref{LDIpDCHScZTgZ}; finally, the last equality follows from the fact that $\mathrm{pr}^{\equiv^{\llbracket 2\rrbracket}}$ is a $\Sigma^{\boldsymbol{\mathcal{A}}^{(2)}}$-homomorphism, according to Remark~\ref{RDVar}.

\textsf{Compatibility with the $0$-target.}

Let $s$ be a sort in $S$ and let $\mathfrak{P}^{(2)}$ be a second-order path in $\mathrm{Pth}_{\boldsymbol{\mathcal{A}}^{(2)},s}$.  

The following equality holds
\allowdisplaybreaks
\begin{multline*}
\mathrm{pr}^{\equiv^{\llbracket 2\rrbracket}}_{s}\left(
\mathrm{ip}^{(2,X)@}_{s}\left(
\mathrm{CH}^{(2)}_{s}\left(
\mathrm{tg}^{0\mathbf{Pth}_{\boldsymbol{\mathcal{A}}^{(2)}}}_{s}\left(
\mathfrak{P}^{(2)}
\right)
\right)\right)\right)
\\=
\mathrm{tg}^{0\mathbf{T}_{\boldsymbol{\mathcal{E}}^{\boldsymbol{\mathcal{A}}^{(2)}}}
(\mathbf{Pth}_{\boldsymbol{\mathcal{A}}^{(2)}})}_{s}\left(
\mathrm{pr}^{\equiv^{\llbracket 2\rrbracket}}_{s}\left(
\mathrm{ip}^{(2,X)@}_{s}\left(
\mathrm{CH}^{(2)}_{s}\left(
\mathfrak{P}^{(2)}
\right)
\right)\right)\right).
\end{multline*}

The proof of this case is similar to that presented for the $0$-source.

\textsf{Compatibility with the $0$-composition.}

Let $s$ be a sort in $S$ and let $\mathfrak{P}^{(2)},\mathfrak{Q}^{(2)}$ be second-order paths in $\mathrm{Pth}_{\boldsymbol{\mathcal{A}}^{(2)},s}$ satisfying  
\[
\mathrm{sc}^{(0,2)}_{s}\left(
\mathfrak{Q}^{(2)}
\right)
=
\mathrm{tg}^{(0,2)}_{s}\left(
\mathfrak{P}^{(2)}
\right).
\]

The following chain of equalities holds
\begin{flushleft}
$\mathrm{pr}^{\equiv^{\llbracket 2\rrbracket}}_{s}\left(
\mathrm{ip}^{(2,X)@}_{s}\left(
\mathrm{CH}^{(2)}_{s}\left(
\mathfrak{Q}^{(2)}
\circ^{0\mathbf{Pth}_{\boldsymbol{\mathcal{A}}^{(2)}}}_{s}
\mathfrak{P}^{(2)}
\right)\right)\right)
$
\allowdisplaybreaks
\begin{align*}
\qquad&=
\mathrm{pr}^{\equiv^{\llbracket 2\rrbracket}}_{s}\left(
\mathrm{ip}^{(2,X)@}_{s}\left(
\mathrm{CH}^{(2)}_{s}\left(
\mathfrak{Q}^{(2)}
\right)\right)
\circ^{0\mathbf{Pth}_{\boldsymbol{\mathcal{A}}^{(2)}}}_{s}
\mathrm{ip}^{(2,X)@}_{s}\left(
\mathrm{CH}^{(2)}_{s}\left(
\mathfrak{P}^{(2)}
\right)\right)
\right)
\tag{1}
\\&=
\mathrm{pr}^{\equiv^{\llbracket 2\rrbracket}}_{s}\left(
\mathrm{ip}^{(2,X)@}_{s}\left(
\mathrm{CH}^{(2)}_{s}\left(
\mathfrak{Q}^{(2)}
\right)
\right)\right)
\circ^{0\mathbf{T}_{\boldsymbol{\mathcal{E}}^{\boldsymbol{\mathcal{A}}^{(2)}}}
(\mathbf{Pth}_{\boldsymbol{\mathcal{A}}^{(2)}})}_{s}
\\&\qquad\qquad\qquad\qquad\qquad\qquad\qquad\qquad
\mathrm{pr}^{\equiv^{\llbracket 2\rrbracket}}_{s}\left(
\mathrm{ip}^{(2,X)@}_{s}\left(
\mathrm{CH}^{(2)}_{s}\left(
\mathfrak{P}^{(2)}
\right)
\right)\right)
.
\tag{2}
\end{align*}
\end{flushleft}

In the just stated chain of equalities, the first equality follows from Lemma~\ref{LDIpDCHCompZ}; finally, the last equality follows from the fact that $\mathrm{pr}^{\equiv^{\llbracket 2\rrbracket}}$ is a $\Sigma^{\boldsymbol{\mathcal{A}}^{(2)}}$-homomorphism, according to Remark~\ref{RDVar}.

\textsf{Compatibility with the second-order rewrite rules in $\mathcal{A}^{(2)}$.}

Let $s$ be a sort in $S$ and let $\mathfrak{p}^{(2)}$ be a second-order rewrite rule in $\mathcal{A}^{(2)}_{s}$.

The following chain of equalities holds.
\allowdisplaybreaks
\begin{align*}
\mathrm{pr}^{\equiv^{\llbracket 2\rrbracket}}_{s}\left(
\mathrm{ip}^{(2,X)@}_{s}\left(
\mathrm{CH}^{(2)}_{s}\left(
\mathfrak{p}^{(2)\mathbf{Pth}_{\boldsymbol{\mathcal{A}}^{(2)}}}
\right)\right)\right)&=
\mathrm{pr}^{\equiv^{\llbracket 2\rrbracket}}_{s}\left(
\mathfrak{p}^{(2)\mathbf{Pth}_{\boldsymbol{\mathcal{A}}^{(2)}}}\right)
\tag{1}
\\&=
\mathfrak{p}^{(2)\mathbf{T}_{\boldsymbol{\mathcal{E}}^{\boldsymbol{\mathcal{A}}^{(2)}}}
(\mathbf{Pth}_{\boldsymbol{\mathcal{A}}^{(2)}})}.
\tag{2}
\end{align*}

In the just stated chain of equalities, the first equality follows from Lemma~\ref{LDIpDCHDEch}; finally, the last equality follows from the fact that $\mathrm{pr}^{\equiv^{\llbracket 2\rrbracket}}$ is a $\Sigma^{\boldsymbol{\mathcal{A}}^{(2)}}$-homomorphism, according to Remark~\ref{RDVar}.

\textsf{Compatibility with the $1$-source.}

Let $s$ be a sort in $S$ and let $\mathfrak{P}^{(2)}$ be a second-order path in $\mathrm{Pth}_{\boldsymbol{\mathcal{A}}^{(2)},s}$. 

The following chain of equalities holds
\begin{flushleft}
$\mathrm{pr}^{\equiv^{\llbracket 2\rrbracket}}_{s}\left(
\mathrm{ip}^{(2,X)@}_{s}\left(
\mathrm{CH}^{(2)}_{s}\left(
\mathrm{sc}^{1\mathbf{Pth}_{\boldsymbol{\mathcal{A}}^{(2)}}}_{s}\left(
\mathfrak{P}^{(2)}
\right)
\right)\right)\right)
$
\allowdisplaybreaks
\begin{align*}
\qquad&=
\mathrm{pr}^{\equiv^{\llbracket 2\rrbracket}}_{s}\left(
\mathrm{sc}^{1\mathbf{Pth}_{\boldsymbol{\mathcal{A}}^{(2)}}}_{s}\left(
\mathrm{ip}^{(2,X)@}_{s}\left(
\mathrm{CH}^{(2)}_{s}\left(
\mathfrak{P}^{(2)}
\right)
\right)\right)\right)
\tag{1}
\\&=
\mathrm{sc}^{1\mathbf{T}_{\boldsymbol{\mathcal{E}}^{\boldsymbol{\mathcal{A}}^{(2)}}}
(\mathbf{Pth}_{\boldsymbol{\mathcal{A}}^{(2)}})}_{s}\left(
\mathrm{pr}^{\equiv^{\llbracket 2\rrbracket}}_{s}\left(
\mathrm{ip}^{(2,X)@}_{s}\left(
\mathrm{CH}^{(2)}_{s}\left(
\mathfrak{P}^{(2)}
\right)
\right)\right)\right).
\tag{2}
\end{align*}
\end{flushleft}

In the just stated chain of equalities, the first equality follows from Lemma~\ref{LDIpDCHScUTgU}; finally, the last equality follows from the fact that $\mathrm{pr}^{\equiv^{\llbracket 2\rrbracket}}$ is a $\Sigma^{\boldsymbol{\mathcal{A}}^{(2)}}$-homomorphism, according to Remark~\ref{RDVar}.

\textsf{Compatibility with the $1$-target.}

Let $s$ be a sort in $S$ and let $\mathfrak{P}^{(2)}$ be a second-order path in $\mathrm{Pth}_{\boldsymbol{\mathcal{A}}^{(2)},s}$.  

The following equality holds
\allowdisplaybreaks
\begin{multline*}
\mathrm{pr}^{\equiv^{\llbracket 2\rrbracket}}_{s}\left(
\mathrm{ip}^{(2,X)@}_{s}\left(
\mathrm{CH}^{(2)}_{s}\left(
\mathrm{tg}^{1\mathbf{Pth}_{\boldsymbol{\mathcal{A}}^{(2)}}}_{s}\left(
\mathfrak{P}^{(2)}
\right)
\right)\right)\right)
\\=
\mathrm{tg}^{1\mathbf{T}_{\boldsymbol{\mathcal{E}}^{\boldsymbol{\mathcal{A}}^{(2)}}}
(\mathbf{Pth}_{\boldsymbol{\mathcal{A}}^{(2)}})}_{s}\left(
\mathrm{pr}^{\equiv^{\llbracket 2\rrbracket}}_{s}\left(
\mathrm{ip}^{(2,X)@}_{s}\left(
\mathrm{CH}^{(2)}_{s}\left(
\mathfrak{P}^{(2)}
\right)
\right)\right)\right).
\end{multline*}

The proof of this case is similar to that presented for the $1$-source.

\textsf{Compatibility with the $1$-composition.}

Let $s$ be a sort in $S$ and let $\mathfrak{P}^{(2)},\mathfrak{Q}^{(2)}$ be second-order paths in $\mathrm{Pth}_{\boldsymbol{\mathcal{A}}^{(2)},s}$ satisfying  
\[
\mathrm{sc}^{([1],2)}_{s}\left(
\mathfrak{Q}^{(2)}
\right)
=
\mathrm{tg}^{([1],2)}_{s}\left(
\mathfrak{P}^{(2)}
\right).
\]

We want to check that 
\allowdisplaybreaks
\begin{multline*}
\mathrm{pr}^{\equiv^{\llbracket 2\rrbracket}}_{s}\left(
\mathrm{ip}^{(2,X)@}_{s}\left(
\mathrm{CH}^{(2)}_{s}\left(
\mathfrak{Q}^{(2)}
\circ^{1\mathbf{Pth}_{\boldsymbol{\mathcal{A}}^{(2)}}}_{s}
\mathfrak{P}^{(2)}
\right)\right)\right)
\\=
\mathrm{pr}^{\equiv^{\llbracket 2\rrbracket}}_{s}\left(
\mathrm{ip}^{(2,X)@}_{s}\left(
\mathrm{CH}^{(2)}_{s}\left(
\mathfrak{Q}^{(2)}
\right)
\right)\right)
\circ^{1\mathbf{T}_{\boldsymbol{\mathcal{E}}^{\boldsymbol{\mathcal{A}}^{(2)}}}
(\mathbf{Pth}_{\boldsymbol{\mathcal{A}}^{(2)}})}_{s}
\\
\mathrm{pr}^{\equiv^{\llbracket 2\rrbracket}}_{s}\left(
\mathrm{ip}^{(2,X)@}_{s}\left(
\mathrm{CH}^{(2)}_{s}\left(
\mathfrak{P}^{(2)}
\right)
\right)\right)
\end{multline*}

Developing both sides of the above equation, and taking into account Proposition~\ref{PDIpDCH}, we have that to check that the following equality between classes holds
\allowdisplaybreaks
\begin{multline*}
\left[
\mathrm{ip}^{(2,X)@}_{s}\left(
\mathrm{CH}^{(2)}_{s}\left(
\mathfrak{Q}^{(2)}
\circ^{1\mathbf{Pth}_{\boldsymbol{\mathcal{A}}^{(2)}}}_{s}
\mathfrak{P}^{(2)}
\right)\right)
\right]_{\equiv^{\llbracket 2\rrbracket}}
\\=
\left[
\mathrm{ip}^{(2,X)@}_{s}\left(
\mathrm{CH}^{(2)}_{s}\left(
\mathfrak{Q}^{(2)}
\right)\right)
\circ^{1\mathbf{Pth}_{\boldsymbol{\mathcal{A}}^{(2)}}}_{s}
\mathrm{ip}^{(2,X)@}_{s}\left(
\mathrm{CH}^{(2)}_{s}\left(
\mathfrak{P}^{(2)}
\right)\right)
\right]_{\equiv^{\llbracket 2\rrbracket}}.
\end{multline*}

But this is equivalent to prove that, for every many-sorted partial $\Sigma^{\boldsymbol{\mathcal{A}}^{(2)}}$-algebra $\mathbf{B}$ in $\mathbf{PAlg}(\boldsymbol{\mathcal{E}}^{\boldsymbol{\mathcal{A}}^{(2)}})$ and every $\Sigma^{\boldsymbol{\mathcal{A}}^{(2)}}$-homomorphism $f\colon\mathbf{Pth}_{\boldsymbol{\mathcal{A}}^{(2)}}\mor \mathbf{B}$ we have that 
\allowdisplaybreaks
\begin{multline*}
f^{\mathrm{Sch}}_{s}\left(
\mathrm{ip}^{(2,X)@}_{s}\left(
\mathrm{CH}^{(2)}_{s}\left(
\mathfrak{Q}^{(2)}
\circ^{1\mathbf{Pth}_{\boldsymbol{\mathcal{A}}^{(2)}}}_{s}
\mathfrak{P}^{(2)}
\right)\right)
\right)
\\=
f^{\mathrm{Sch}}_{s}\left(
\mathrm{ip}^{(2,X)@}_{s}\left(
\mathrm{CH}^{(2)}_{s}\left(
\mathfrak{Q}^{(2)}
\right)\right)
\circ^{1\mathbf{Pth}_{\boldsymbol{\mathcal{A}}^{(2)}}}_{s}
\mathrm{ip}^{(2,X)@}_{s}\left(
\mathrm{CH}^{(2)}_{s}\left(
\mathfrak{P}^{(2)}
\right)\right)
\right).
\end{multline*}

Let us recall, by Proposition~\ref{PDIpDCH}, that we are dealing with  proper second-order paths. Hence, the last mentioned equality reduces to
\allowdisplaybreaks
\begin{multline*}
f_{s}\left(
\mathrm{ip}^{(2,X)@}_{s}\left(
\mathrm{CH}^{(2)}_{s}\left(
\mathfrak{Q}^{(2)}
\circ^{1\mathbf{Pth}_{\boldsymbol{\mathcal{A}}^{(2)}}}_{s}
\mathfrak{P}^{(2)}
\right)\right)
\right)
\\=
f_{s}\left(
\mathrm{ip}^{(2,X)@}_{s}\left(
\mathrm{CH}^{(2)}_{s}\left(
\mathfrak{Q}^{(2)}
\right)\right)
\circ^{1\mathbf{Pth}_{\boldsymbol{\mathcal{A}}^{(2)}}}_{s}
\mathrm{ip}^{(2,X)@}_{s}\left(
\mathrm{CH}^{(2)}_{s}\left(
\mathfrak{P}^{(2)}
\right)\right)
\right).
\tag{E2}\label{PDVarKerE2}
\end{multline*}

Before presenting a proof for the general case, we will see a proof of Equation~\ref{PDVarKerE2} when one of the second-order paths involved in the $1$-composition is a $(2,[1])$-identity second-order path.

\begin{claim}\label{CDVarKer} Let $\mathfrak{P}^{(2)}$ and $\mathfrak{Q}^{(2)}$ be two second-order paths in $\mathrm{Pth}_{\boldsymbol{\mathcal{A}}^{(2)},s}$ satisfying that $\mathrm{sc}^{([1],2)}_{s}(\mathfrak{Q}^{(2)})=\mathrm{tg}^{([1],2)}_{s}(\mathfrak{P}^{(2)})$. If either $\mathfrak{P}^{(2)}$ or $\mathfrak{Q}^{(2)}$ is a $(2,[1])$-identity second-order path, then for every many-sorted partial $\Sigma^{\boldsymbol{\mathcal{A}}^{(2)}}$-algebra $\mathbf{B}$ in $\mathbf{PAlg}(\boldsymbol{\mathcal{E}}^{\boldsymbol{\mathcal{A}}^{(2)}})$ and every $\Sigma^{\boldsymbol{\mathcal{A}}^{(2)}}$-homomorphism $f\colon \mathbf{Pth}_{\boldsymbol{\mathcal{A}}^{(2)}}\mor \mathbf{B}$, Equation~\ref{PDVarKerE2} holds.
\end{claim}

Assume without loss of generality that $\mathfrak{P}^{(2)}$ is a $(2,[1])$-identity second-order path, then the following chain of equalities holds.
\begin{flushleft}
$f_{s}\left(
\mathrm{ip}^{(2,X)@}_{s}\left(
\mathrm{CH}^{(2)}_{s}\left(
\mathfrak{Q}^{(2)}
\circ^{1\mathbf{Pth}_{\boldsymbol{\mathcal{A}}^{(2)}}}_{s}
\mathfrak{P}^{(2)}
\right)\right)
\right)$
\allowdisplaybreaks
\begin{align*}
&=f_{s}\left(
\mathrm{ip}^{(2,X)@}_{s}\left(
\mathrm{CH}^{(2)}_{s}\left(
\mathfrak{Q}^{(2)}
\right)\right)
\right)
\tag{1}
\\&=f_{s}\left(
\mathrm{ip}^{(2,X)@}_{s}\left(
\mathrm{CH}^{(2)}_{s}\left(
\mathfrak{Q}^{(2)}
\right)\right)
\circ^{1\mathbf{Pth}_{\boldsymbol{\mathcal{A}}^{(2)}}}_{s}
\mathfrak{P}^{(2)}
\right)
\tag{2}
\\&=f_{s}\left(
\mathrm{ip}^{(2,X)@}_{s}\left(
\mathrm{CH}^{(2)}_{s}\left(
\mathfrak{Q}^{(2)}
\right)\right)
\right)
\circ^{1\mathbf{B}}_{s}
f_{s}\left(
\mathfrak{P}^{(2)}
\right)
\tag{3}
\\&=
f_{s}\left(
\mathrm{ip}^{(2,X)@}_{s}\left(
\mathrm{CH}^{(2)}_{s}\left(
\mathfrak{Q}^{(2)}
\right)\right)
\right)
\circ^{1\mathbf{B}}_{s}
f_{s}\left(
\mathrm{ip}^{(2,X)@}_{s}\left(
\mathrm{CH}^{(2)}_{s}\left(
\mathfrak{P}^{(2)}
\right)\right)
\right).
\tag{4}
\end{align*}
\end{flushleft}

Note that since $\mathfrak{P}^{(2)}$ is a $(2,[1])$-identity second-order path, this second-order path must be equal to the $(2,[1])$-identity second-order path on the $([1],2)$-source of $\mathfrak{Q}^{(2)}$, so the $1$-composition of the second-order paths is well-defined. This justifies the first equality. The second equality follows from the fact that, according to Proposition~\ref{PDIpDCH}, $\mathrm{ip}^{(2,X)@}_{s}(\mathrm{CH}^{(2)}_{s}(\mathfrak{Q}^{(2)}))$ is a second-order path in $[\mathfrak{Q}^{(2)}]_{s}$. Thus, by Lemma~\ref{LDCH}, the second-order paths $\mathrm{ip}^{(2,X)@}_{s}(\mathrm{CH}^{(2)}_{s}(\mathfrak{Q}^{(2)}))$ and $\mathfrak{Q}^{(2)}$ have the same $([1],2)$-source. Moreover, since $\mathfrak{P}^{(2)}$ is the  $(2,[1])$-identity second-order path on the $([1],2)$-source of $\mathfrak{Q}^{(2)}$, we have that
\[
\mathrm{ip}^{(2,X)@}_{s}\left(
\mathrm{CH}^{(2)}_{s}\left(
\mathfrak{Q}^{(2)}
\right)\right)
=
\mathrm{ip}^{(2,X)@}_{s}\left(
\mathrm{CH}^{(2)}_{s}\left(
\mathfrak{Q}^{(2)}
\right)\right)
\circ^{1\mathbf{Pth}_{\boldsymbol{\mathcal{A}}^{(2)}}}_{s}
\mathfrak{P}^{(2)};
\]
The third equality follows from the fact that, by assumption, $f$ is a $\Sigma^{\boldsymbol{\mathcal{A}}^{(2)}}$-homomorphism; finally, the last equality follows from Corollary~\ref{CDIpDCHOneStep}. In this regard, since $\mathfrak{P}^{(2)}$ is a $(2,[1])$-identity  second-order path, we have that $\mathrm{ip}^{(2,X)@}_{s}(
\mathrm{CH}^{(2)}_{s}(
\mathfrak{P}^{(2)}
))=\mathfrak{P}^{(2)}$.

The same argument applies in case $\mathfrak{Q}^{(2)}$ is a $(2,[1])$-identity  second-order path. In this case, we will argue taking into account that $\mathfrak{Q}^{(2)}$ must be the $(2,[1])$-identity  second-order path on the $([1],2)$-target of $\mathfrak{P}^{(2)}$.

This proves Claim~\ref{CDVarKer}.

We are now in position to prove the general case on $\mathfrak{Q}^{(2)}\circ^{1\mathbf{Pth}_{\boldsymbol{\mathcal{A}}^{(2)}}}_{s}\mathfrak{P}^{(2)}$ by Artinian induction on $(\coprod\mathrm{Pth}_{\boldsymbol{\mathcal{A}}^{(2)}}, \leq_{\mathbf{Pth}_{\boldsymbol{\mathcal{A}}^{(2)}}})$.

\textsf{Base step of the Artinian induction.}

Let $(\mathfrak{Q}^{(2)}\circ^{1\mathbf{Pth}_{\boldsymbol{\mathcal{A}}^{(2)}}}_{s}\mathfrak{P}^{(2)}, s)$ be a minimal element in $(\coprod\mathrm{Pth}_{\boldsymbol{\mathcal{A}}^{(2)}}, \leq_{\mathbf{Pth}_{\boldsymbol{\mathcal{A}}^{(2)}}})$. Then by Proposition~\ref{PDMinimal}, the second-order path $\mathfrak{Q}^{(2)}\circ^{1\mathbf{Pth}_{\boldsymbol{\mathcal{A}}^{(2)}}}_{s}\mathfrak{P}^{(2)}$ is either a $(2,[1])$-identity second-order path or a  second-order echelon. In any case, either $\mathfrak{P}^{(2)}$ or $\mathfrak{Q}^{(2)}$ must be a $(2,[1])$-identity second-order path. The statement follows by Claim~\ref{CDVarKer}.

\textsf{Inductive step of the Artinian induction.}

Let $(\mathfrak{Q}^{(2)}\circ^{1\mathbf{Pth}_{\boldsymbol{\mathcal{A}}^{(2)}}}_{s}\mathfrak{P}^{(2)}, s)$ be a non-minimal element in $(\coprod\mathrm{Pth}_{\boldsymbol{\mathcal{A}}^{(2)}}, \leq_{\mathbf{Pth}_{\boldsymbol{\mathcal{A}}^{(2)}}})$. Let us suppose that, for every sort $t\in S$ and every second-order path $\mathfrak{Q}'^{(2)}\circ^{1\mathbf{Pth}_{\boldsymbol{\mathcal{A}}^{(2)}}}_{t}\mathfrak{P}'^{(2)}$ in $\mathrm{Pth}_{\boldsymbol{\mathcal{A}}^{(2)},t}$, if $(\mathfrak{Q}'^{(2)}\circ^{1\mathbf{Pth}_{\boldsymbol{\mathcal{A}}^{(2)}}}_{t}\mathfrak{P}'^{(2)},t)<_{\mathbf{Pth}_{\boldsymbol{\mathcal{A}}^{(2)}}} (\mathfrak{Q}^{(2)}\circ^{1\mathbf{Pth}_{\boldsymbol{\mathcal{A}}^{(2)}}}_{s}\mathfrak{P}^{(2)}, s)$, then the statement holds for $\mathfrak{Q}'^{(2)}\circ^{1\mathbf{Pth}_{\boldsymbol{\mathcal{A}}^{(2)}}}_{t}\mathfrak{P}'^{(2)}$, i.e., for every many-sorted partial $\Sigma^{\boldsymbol{\mathcal{A}}^{(2)}}$-algebra $\mathbf{B}$ in $\mathbf{PAlg}(\boldsymbol{\mathcal{E}}^{\boldsymbol{\mathcal{A}}^{(2)}})$ and every $\Sigma^{\boldsymbol{\mathcal{A}}^{(2)}}$-homomorphism $f\colon \mathbf{Pth}_{\boldsymbol{\mathcal{A}}^{(2)}}\mor \mathbf{B}$, the following equality holds
\allowdisplaybreaks
\begin{multline*}
f_{t}\left(
\mathrm{ip}^{(2,X)@}_{t}\left(
\mathrm{CH}^{(2)}_{t}\left(
\mathfrak{Q}'^{(2)}
\circ^{1\mathbf{Pth}_{\boldsymbol{\mathcal{A}}^{(2)}}}_{t}
\mathfrak{P}'^{(2)}
\right)\right)
\right)
\\=
f_{t}\left(
\mathrm{ip}^{(2,X)@}_{t}\left(
\mathrm{CH}^{(2)}_{t}\left(
\mathfrak{Q}'^{(2)}
\right)\right)
\circ^{1\mathbf{Pth}_{\boldsymbol{\mathcal{A}}^{(2)}}}_{t}
\mathrm{ip}^{(2,X)@}_{t}\left(
\mathrm{CH}^{(2)}_{t}\left(
\mathfrak{P}'^{(2)}
\right)\right)
\right).
\end{multline*}

Since $(\mathfrak{Q}^{(2)}\circ^{1\mathbf{Pth}_{\boldsymbol{\mathcal{A}}^{(2)}}}_{s}\mathfrak{P}^{(2)}, s)$ is a non-minimal element in $(\coprod\mathrm{Pth}_{\boldsymbol{\mathcal{A}}^{(2)}}, \leq_{\mathbf{Pth}_{\boldsymbol{\mathcal{A}}^{(2)}}})$ and, by Claim~\ref{CDVarKer}, we can assume that neither $\mathfrak{Q}^{(2)}$ nor $\mathfrak{P}^{(2)}$ are $(2,[1])$-identity second-order paths, we have, by Lemma~\ref{LDOrdI} that $\mathfrak{Q}^{(2)}\circ^{1\mathbf{Pth}_{\boldsymbol{\mathcal{A}}^{(2)}}}_{s}\mathfrak{P}^{(2)}$ is either (1) a second-order path of length strictly greater than one containing at least one  second-order echelon or (2) an echelonless second-order path.

If~(1), then let $i\in\bb{\mathfrak{Q}^{(2)}\circ^{1\mathbf{Pth}_{\boldsymbol{\mathcal{A}}^{(2)}}}_{s}\mathfrak{P}^{(2)}}$ be the first index for which the one-step subpath $(\mathfrak{Q}^{(2)}\circ^{1\mathbf{Pth}_{\boldsymbol{\mathcal{A}}^{(2)}}}_{s}\mathfrak{P}^{(2)})^{i,i}$ of $\mathfrak{Q}^{(2)}\circ^{1\mathbf{Pth}_{\boldsymbol{\mathcal{A}}^{(2)}}}_{s}\mathfrak{P}^{(2)}$ is a  second-order echelon. We distinguish the cases~(1.1) $i=0$ and (1.2) $i\neq 0$.

If~(1.1), i.e., if $i=0$, since we are assuming that $\mathfrak{P}^{(2)}$ is not a $(2,[1])$-identity second-order path, we have that $\mathfrak{P}^{(2)}$ has a  second-order echelon on its first step. Then it could be the case that either (1.1.1) $\mathfrak{P}^{(2)}$
 is a  second-order echelon, or (1.1.2) $\mathfrak{P}^{(2)}$ is a second-order path of length strictly greater than one containing a  second-order echelon on its first step.
 
If~(1.1.1), i.e., if we consider the case in which $\mathfrak{P}^{(2)}$ is a  second-order echelon, then the value of the second-order Curry-Howard mapping at $\mathfrak{Q}^{(2)}
\circ^{1\mathbf{Pth}_{\boldsymbol{\mathcal{A}}^{(2)}}}_{s}
\mathfrak{P}^{(2)}$ is given by
\[
\mathrm{CH}^{(2)}_{s}\left(
\mathfrak{Q}^{(2)}
\circ^{1\mathbf{Pth}_{\boldsymbol{\mathcal{A}}^{(2)}}}_{s}
\mathfrak{P}^{(2)}
\right)
=
\mathrm{CH}^{(2)}_{s}\left(
\mathfrak{Q}^{(2)}
\right)
\circ^{1\mathbf{T}_{\Sigma^{\boldsymbol{\mathcal{A}}^{(2)}}}(X)}_{s}
\mathrm{CH}^{(2)}_{s}\left(
\mathfrak{P}^{(2)}
\right).
\]

Therefore, the following chain of equalities holds
\begin{flushleft}
$f_{s}\left(
\mathrm{ip}^{(2,X)@}_{s}\left(
\mathrm{CH}^{(2)}_{s}\left(
\mathfrak{Q}^{(2)}
\circ^{1\mathbf{Pth}_{\boldsymbol{\mathcal{A}}^{(2)}}}_{s}
\mathfrak{P}^{(2)}
\right)\right)
\right)$
\allowdisplaybreaks
\begin{align*}
&=
f_{s}\left(
\mathrm{ip}^{(2,X)@}_{s}\left(
\mathrm{CH}^{(2)}_{s}\left(
\mathfrak{Q}^{(2)}\right)
\circ^{1\mathbf{T}_{\Sigma^{\boldsymbol{\mathcal{A}}^{(2)}}}(X)}_{s}
\mathrm{CH}^{(2)}_{s}\left(
\mathfrak{P}^{(2)}
\right)\right)
\right)
\tag{1}
\\&=
f_{s}\left(
\mathrm{ip}^{(2,X)@}_{s}\left(
\mathrm{CH}^{(2)}_{s}\left(
\mathfrak{Q}^{(2)}\right)\right)
\circ^{1\mathbf{F}_{\Sigma^{\boldsymbol{\mathcal{A}}^{(2)}}}(\mathbf{Pth}_{\boldsymbol{\mathcal{A}}^{(2)}})}_{s}
\mathrm{ip}^{(2,X)@}_{s}\left(
\mathrm{CH}^{(2)}_{s}\left(
\mathfrak{P}^{(2)}
\right)\right)
\right)
\tag{2}
\\&=
f_{s}\left(
\mathrm{ip}^{(2,X)@}_{s}\left(
\mathrm{CH}^{(2)}_{s}\left(
\mathfrak{Q}^{(2)}\right)\right)
\circ^{1\mathbf{Pth}_{\boldsymbol{\mathcal{A}}^{(2)}}}_{s}
\mathrm{ip}^{(2,X)@}_{s}\left(
\mathrm{CH}^{(2)}_{s}\left(
\mathfrak{P}^{(2)}
\right)\right)
\right)
\tag{3}
\\&=
f_{s}\left(
\mathrm{ip}^{(2,X)@}_{s}\left(
\mathrm{CH}^{(2)}_{s}\left(
\mathfrak{Q}^{(2)}\right)\right)\right)
\circ^{1\mathbf{B}}_{s}
f_{s}\left(
\mathrm{ip}^{(2,X)@}_{s}\left(
\mathrm{CH}^{(2)}_{s}\left(
\mathfrak{P}^{(2)}
\right)\right)
\right).
\tag{4}
\end{align*}
\end{flushleft}

In the just stated chain of equalities, the first equality follows from previous discussion on the value of the second-order Curry-Howard mapping at $\mathfrak{Q}^{(2)}
\circ^{1\mathbf{Pth}_{\boldsymbol{\mathcal{A}}^{(2)}}}_{s}
\mathfrak{P}^{(2)}$; the second equality follow from the fact that $\mathrm{ip}^{(2,X)@}$ is a $\Sigma^{\boldsymbol{\mathcal{A}}^{(2)}}$-homomorphism, according to Definition~\ref{DDIp}; the third equality follows from the fact that, by Proposition~\ref{PDIpDCH}, both $\mathrm{ip}^{(2,X)@}_{s}(
\mathrm{CH}^{(2)}_{s}(
\mathfrak{Q}^{(2)}
))$ and $\mathrm{ip}^{(2,X)@}_{s}(
\mathrm{CH}^{(2)}_{s}(
\mathfrak{P}^{(2)}
))$ are second-order paths. Thus, the interpretation of the operation symbol for the $1$-composition $\circ^{1}_{s}$ in the $\Sigma^{\boldsymbol{\mathcal{A}}^{(2)}}$-algebra $\mathbf{F}_{\Sigma^{\boldsymbol{\mathcal{A}}^{(2)}}}(\mathbf{Pth}_{\boldsymbol{\mathcal{A}}^{(2)}})$ becomes the corresponding interpretation of the $1$-composition $\circ^{1}_{s}$ in the many-sorted partial $\Sigma^{\boldsymbol{\mathcal{A}}^{(2)}}$-algebra $\mathbf{Pth}_{\boldsymbol{\mathcal{A}}^{(2)}}$; finally, the last equality follows from the fact that, by assumption, $f$ is a $\Sigma^{\boldsymbol{\mathcal{A}}^{(2)}}$-homomorphism.

The case~(1.1.1) of $\mathfrak{P}^{(2)}$ being a  second-order echelon follows.

If~(1.1.2), i.e., if we are in the case in which $\mathfrak{P}^{(2)}$  is a second-order path of length strictly greater than one containing a  second-order echelon on its first step, then the value of the second-order Curry-Howard mapping at $\mathfrak{P}^{(2)}$ is given by
\[
\mathrm{CH}^{(2)}_{s}\left(
\mathfrak{P}^{(2)}
\right)
=
\mathrm{CH}^{(2)}_{s}\left(
\mathfrak{P}^{(2),1,\bb{\mathfrak{P}^{(2)}}-1}
\right)
\circ^{1\mathbf{T}_{\Sigma^{\boldsymbol{\mathcal{A}}^{(2)}}}(X)}_{s}
\mathrm{CH}^{(2)}_{s}\left(
\mathfrak{P}^{(2),0,0}
\right).
\]

Moreover, $(\mathfrak{Q}^{(2)}
\circ^{1\mathbf{Pth}_{\boldsymbol{\mathcal{A}}^{(2)}}}_{s}
\mathfrak{P}^{(2)})^{1,\bb{\mathfrak{Q}^{(2)}
\circ^{1\mathbf{Pth}_{\boldsymbol{\mathcal{A}}^{(2)}}}_{s}
\mathfrak{P}^{(2)}}-1}=\mathfrak{Q}^{(2)}\circ^{1\mathbf{Pth}_{\boldsymbol{\mathcal{A}}^{(2)}}}_{s}\mathfrak{P}^{(2),1,\bb{\mathfrak{P}^{(2))}-1}}$. Hence, the value of the second-order Curry-Howard mapping at $\mathfrak{Q}^{(2)}
\circ^{1\mathbf{Pth}_{\boldsymbol{\mathcal{A}}^{(2)}}}_{s}
\mathfrak{P}^{(2)}$ is given by 
\allowdisplaybreaks
\begin{multline*}
\mathrm{CH}^{(2)}_{s}\left(
\mathfrak{Q}^{(2)}
\circ^{1\mathbf{Pth}_{\boldsymbol{\mathcal{A}}^{(2)}}}_{s}
\mathfrak{P}^{(2)}
\right)
\\=
\mathrm{CH}^{(2)}_{s}\left(
\mathfrak{Q}^{(2)}\circ^{1\mathbf{Pth}_{\boldsymbol{\mathcal{A}}^{(2)}}}_{s}\mathfrak{P}^{(2),1,\bb{\mathfrak{P}^{(2))}-1}}
\right)
\circ^{1\mathbf{T}_{\Sigma^{\boldsymbol{\mathcal{A}}^{(2)}}}(X)}_{s}
\mathrm{CH}^{(2)}_{s}\left(
\mathfrak{P}^{(2),0,0}
\right).
\end{multline*}

Therefore, the following chain of equalities holds
\begin{flushleft}
$f_{s}\left(
\mathrm{ip}^{(2,X)@}_{s}\left(
\mathrm{CH}^{(2)}_{s}\left(
\mathfrak{Q}^{(2)}
\circ^{1\mathbf{Pth}_{\boldsymbol{\mathcal{A}}^{(2)}}}_{s}
\mathfrak{P}^{(2)}
\right)\right)
\right)$
\allowdisplaybreaks
\begin{align*}
&=f_{s}\left(
\mathrm{ip}^{(2,X)@}_{s}\left(
\mathrm{CH}^{(2)}_{s}\left(
\mathfrak{Q}^{(2)}\circ^{1\mathbf{Pth}_{\boldsymbol{\mathcal{A}}^{(2)}}}_{s}\mathfrak{P}^{(2),1,\bb{\mathfrak{P}^{(2))}-1}}
\right)
\right.\right.
\\&\qquad\qquad\qquad\qquad\qquad\qquad\qquad\qquad
\left.\left.
\circ^{1\mathbf{T}_{\Sigma^{\boldsymbol{\mathcal{A}}^{(2)}}}(X)}_{s}
\mathrm{CH}^{(2)}_{s}\left(
\mathfrak{P}^{(2),0,0}
\right)
\right)
\right)
\tag{1}
\\&=f_{s}\left(
\mathrm{ip}^{(2,X)@}_{s}\left(
\mathrm{CH}^{(2)}_{s}\left(
\mathfrak{Q}^{(2)}\circ^{1\mathbf{Pth}_{\boldsymbol{\mathcal{A}}^{(2)}}}_{s}\mathfrak{P}^{(2),1,\bb{\mathfrak{P}^{(2))}-1}}
\right)
\right)\right.
\\&\qquad\qquad\qquad\qquad\qquad\qquad
\left.
\circ^{1\mathbf{F}_{\Sigma^{\boldsymbol{\mathcal{A}}^{(2)}}}(\mathbf{Pth}_{\boldsymbol{\mathcal{A}}^{(2)}})}_{s}
\mathrm{ip}^{(2,X)@}_{s}\left(
\mathrm{CH}^{(2)}_{s}\left(
\mathfrak{P}^{(2),0,0}
\right)
\right)
\right)
\tag{2}
\\&=f_{s}\left(
\mathrm{ip}^{(2,X)@}_{s}\left(
\mathrm{CH}^{(2)}_{s}\left(
\mathfrak{Q}^{(2)}\circ^{1\mathbf{Pth}_{\boldsymbol{\mathcal{A}}^{(2)}}}_{s}\mathfrak{P}^{(2),1,\bb{\mathfrak{P}^{(2))}-1}}
\right)
\right)\right.
\\&\qquad\qquad\qquad\qquad\qquad\qquad\qquad\qquad
\left.
\circ^{1\mathbf{Pth}_{\boldsymbol{\mathcal{A}}^{(2)}}}_{s}
\mathrm{ip}^{(2,X)@}_{s}\left(
\mathrm{CH}^{(2)}_{s}\left(
\mathfrak{P}^{(2),0,0}
\right)
\right)
\right)
\tag{3}
\\&=f_{s}\left(
\mathrm{ip}^{(2,X)@}_{s}\left(
\mathrm{CH}^{(2)}_{s}\left(
\mathfrak{Q}^{(2)}\circ^{1\mathbf{Pth}_{\boldsymbol{\mathcal{A}}^{(2)}}}_{s}\mathfrak{P}^{(2),1,\bb{\mathfrak{P}^{(2))}-1}}
\right)
\right)\right)
\\&\qquad\qquad\qquad\qquad\qquad\qquad\qquad\qquad
\circ^{1\mathbf{B}}_{s}
f_{s}\left(
\mathrm{ip}^{(2,X)@}_{s}\left(
\mathrm{CH}^{(2)}_{s}\left(
\mathfrak{P}^{(2),0,0}
\right)
\right)
\right)
\tag{4}
\\&=\left(f_{s}\left(
\mathrm{ip}^{(2,X)@}_{s}\left(
\mathrm{CH}^{(2)}_{s}\left(
\mathfrak{Q}^{(2)}
\right)\right)\right)
\right.
\\&\qquad\qquad\qquad\quad
\left.
\circ^{1\mathbf{B}}_{s}
f_{s}\left(
\mathrm{ip}^{(2,X)@}_{s}\left(
\mathrm{CH}^{(2)}_{s}\left(
\mathfrak{P}^{(2),1,\bb{\mathfrak{P}^{(2))}-1}}
\right)
\right)\right)\right)
\\&\qquad\qquad\qquad\qquad\qquad\qquad\qquad\qquad
\circ^{1\mathbf{B}}_{s}
f_{s}\left(
\mathrm{ip}^{(2,X)@}_{s}\left(
\mathrm{CH}^{(2)}_{s}\left(
\mathfrak{P}^{(2),0,0}
\right)
\right)
\right)
\tag{5}
\\&=f_{s}\left(
\mathrm{ip}^{(2,X)@}_{s}\left(
\mathrm{CH}^{(2)}_{s}\left(
\mathfrak{Q}^{(2)}
\right)\right)\right)
\circ^{1\mathbf{B}}_{s}
\\&\qquad\qquad\qquad\quad
\left(
f_{s}\left(
\mathrm{ip}^{(2,X)@}_{s}\left(
\mathrm{CH}^{(2)}_{s}\left(
\mathfrak{P}^{(2),1,\bb{\mathfrak{P}^{(2))}-1}}
\right)
\right)\right)
\right.
\\&\qquad\qquad\qquad\qquad\qquad\qquad\qquad\qquad
\circ^{1\mathbf{B}}_{s}
\left.
f_{s}\left(
\mathrm{ip}^{(2,X)@}_{s}\left(
\mathrm{CH}^{(2)}_{s}\left(
\mathfrak{P}^{(2),0,0}
\right)
\right)
\right)\right)
\tag{6}
\\&=f_{s}\left(
\mathrm{ip}^{(2,X)@}_{s}\left(
\mathrm{CH}^{(2)}_{s}\left(
\mathfrak{Q}^{(2)}
\right)\right)\right)
\circ^{1\mathbf{B}}_{s}
\\&\qquad\qquad\qquad\quad
\left(
f_{s}\left(
\mathrm{ip}^{(2,X)@}_{s}\left(
\mathrm{CH}^{(2)}_{s}\left(
\mathfrak{P}^{(2),1,\bb{\mathfrak{P}^{(2))}-1}}
\right)
\right)\right)
\right.
\\&\qquad\qquad\qquad\qquad\qquad\qquad\qquad
\circ^{1\mathbf{\mathbf{Pth}_{\boldsymbol{\mathcal{A}}^{(2)}}}}_{s}
\left.\left.
\mathrm{ip}^{(2,X)@}_{s}\left(
\mathrm{CH}^{(2)}_{s}\left(
\mathfrak{P}^{(2),0,0}
\right)\right)
\right)\right)
\tag{7}
\\&=f_{s}\left(
\mathrm{ip}^{(2,X)@}_{s}\left(
\mathrm{CH}^{(2)}_{s}\left(
\mathfrak{Q}^{(2)}
\right)\right)\right)
\circ^{1\mathbf{B}}_{s}
\\&\qquad\qquad\qquad\quad
\left(
f_{s}\left(
\mathrm{ip}^{(2,X)@}_{s}\left(
\mathrm{CH}^{(2)}_{s}\left(
\mathfrak{P}^{(2),1,\bb{\mathfrak{P}^{(2))}-1}}
\right)
\right)\right)
\right.
\\&\qquad\qquad\qquad\qquad\qquad\qquad
\circ^{1\mathbf{F}_{\Sigma^{\boldsymbol{\mathcal{A}}^{(2)}}}(\mathbf{\mathbf{Pth}_{\boldsymbol{\mathcal{A}}^{(2)}}})}_{s}
\left.\left.
\mathrm{ip}^{(2,X)@}_{s}\left(
\mathrm{CH}^{(2)}_{s}\left(
\mathfrak{P}^{(2),0,0}
\right)\right)
\right)\right)
\tag{8}
\\&=f_{s}\left(
\mathrm{ip}^{(2,X)@}_{s}\left(
\mathrm{CH}^{(2)}_{s}\left(
\mathfrak{Q}^{(2)}
\right)\right)\right)
\circ^{1\mathbf{B}}_{s}
\\&\qquad\qquad\qquad\qquad\qquad
\left(
f_{s}\left(
\mathrm{ip}^{(2,X)@}_{s}\left(
\mathrm{CH}^{(2)}_{s}\left(
\mathfrak{P}^{(2),1,\bb{\mathfrak{P}^{(2))}-1}}
\right)
\right.\right.\right.
\\&\qquad\qquad\qquad\qquad\qquad\qquad\qquad\qquad\qquad
\left.\left.\left.
\circ^{1\mathbf{T}_{\Sigma^{\boldsymbol{\mathcal{A}}^{(2)}}}(X)}_{s}
\mathrm{CH}^{(2)}_{s}\left(
\mathfrak{P}^{(2),0,0}
\right)
\right)
\right)
\right)
\tag{9}
\\&=
f_{s}\left(
\mathrm{ip}^{(2,X)@}_{s}\left(
\mathrm{CH}^{(2)}_{s}\left(
\mathfrak{Q}^{(2)}\right)\right)\right)
\circ^{1\mathbf{B}}_{s}
f_{s}\left(
\mathrm{ip}^{(2,X)@}_{s}\left(
\mathrm{CH}^{(2)}_{s}\left(
\mathfrak{P}^{(2)}
\right)\right)
\right).
\tag{10}
\end{align*}
\end{flushleft}

In the just stated chain of equalities, the first equality follows from the previous discussion on the value of the second-order Curry-Howard mapping at $\mathfrak{Q}^{(2)}\circ^{1\mathbf{Pth}_{\boldsymbol{\mathcal{A}}^{(2)}}}_{s}\mathfrak{P}^{(2)}$; the second equality follow from the fact that $\mathrm{ip}^{(2,X)@}$ is a $\Sigma^{\boldsymbol{\mathcal{A}}^{(2)}}$-homomorphism, according to Definition~\ref{DDIp}; the third equality follows from the fact that, by Proposition~\ref{PDIpDCH}, both $\mathrm{ip}^{(2,X)@}_{s}(
\mathrm{CH}^{(2)}_{s}(
\mathfrak{Q}^{(2)}\circ^{1\mathbf{Pth}_{\boldsymbol{\mathcal{A}}^{(2)}}}_{s}
\mathfrak{P}^{(2),1,\bb{\mathfrak{P}^{(2)}}-1}
))$ and $\mathrm{ip}^{(2,X)@}_{s}(
\mathrm{CH}^{(2)}_{s}(
\mathfrak{P}^{(2),0,0}
))$ are second-order paths. Thus, the interpretation of the operation symbol for the $1$-composition $\circ^{1}_{s}$ in the $\Sigma^{\boldsymbol{\mathcal{A}}^{(2)}}$-algebra $\mathbf{F}_{\Sigma^{\boldsymbol{\mathcal{A}}^{(2)}}}(\mathbf{Pth}_{\boldsymbol{\mathcal{A}}^{(2)}})$ becomes the corresponding interpretation of the $1$-composition $\circ^{1}_{s}$ in the many-sorted partial $\Sigma^{\boldsymbol{\mathcal{A}}^{(2)}}$-algebra $\mathbf{Pth}_{\boldsymbol{\mathcal{A}}^{(2)}}$; the fourth equality follows from the fact that, by assumption $f$ is a $\Sigma^{\boldsymbol{\mathcal{A}}^{(2)}}$-homomorphism; the fifth equality follows by induction. Let us note that the pair $(\mathfrak{Q}^{(2)}\circ^{1\mathbf{Pth}_{\boldsymbol{\mathcal{A}}^{(2)}}}_{s}
\mathfrak{P}^{(2),1,\bb{\mathfrak{P}^{(2)}}-1},s)$ $\prec_{\mathbf{Pth}_{\boldsymbol{\mathcal{A}}^{(2)}}}$-precedes $(\mathfrak{Q}^{(2)}\circ^{1\mathbf{Pth}_{\boldsymbol{\mathcal{A}}^{(2)}}}_{s}\mathfrak{P}^{(2)},s)$. Hence, we have that 
\begin{multline*}
f_{s}\left(
\mathrm{ip}^{(2,X)@}_{s}\left(
\mathrm{CH}^{(2)}_{s}\left(
\mathfrak{Q}^{(2)}
\circ^{1\mathbf{Pth}_{\boldsymbol{\mathcal{A}}^{(2)}}}_{s}
\mathfrak{P}^{(2),1,\bb{\mathfrak{P}^{(2)}}-1}
\right)\right)
\right)
\\=
f_{s}\left(
\mathrm{ip}^{(2,X)@}_{s}\left(
\mathrm{CH}^{(2)}_{s}\left(
\mathfrak{Q}^{(2)}
\right)\right)\right)
\circ^{1\mathbf{B}}_{s}
f_{s}\left(
\mathrm{ip}^{(2,X)@}_{s}\left(
\mathrm{CH}^{(2)}_{s}\left(
\mathfrak{P}^{(2),1,\bb{\mathfrak{P}^{(2)}}-1}
\right)\right)
\right);
\end{multline*}
the sixth equality follows from the fact that, since $\mathbf{B}$ is a partial many-sorted $\Sigma^{\boldsymbol{\mathcal{A}}^{(2)}}$-algebra in $\mathbf{PAlg}(\boldsymbol{\mathcal{E}}^{\boldsymbol{\mathcal{A}}^{(2)}})$, the $1$-composition is associative; the seventh equality follows from the fact that, by assumption, $f$ is a $\Sigma^{\boldsymbol{\mathcal{A}}^{(2)}}$-homomorphism; the eighth equality follows from the fact that,  by Proposition~\ref{PDIpDCH}, both $\mathrm{ip}^{(2,X)@}_{s}(
\mathrm{CH}^{(2)}_{s}(
\mathfrak{P}^{(2),1,\bb{\mathfrak{P}^{(2)}}-1}
))$ and $\mathrm{ip}^{(2,X)@}_{s}(
\mathrm{CH}^{(2)}_{s}(
\mathfrak{P}^{(2),0,0}
))$ are second-order paths. Thus, the interpretation of the operation symbol for the $1$-composition $\circ^{1}_{s}$ in the $\Sigma^{\boldsymbol{\mathcal{A}}^{(2)}}$-algebra $\mathbf{F}_{\Sigma^{\boldsymbol{\mathcal{A}}^{(2)}}}(\mathbf{Pth}_{\boldsymbol{\mathcal{A}}^{(2)}})$ becomes the corresponding interpretation of the $1$-composition $\circ^{1}_{s}$ in the many-sorted partial $\Sigma^{\boldsymbol{\mathcal{A}}^{(2)}}$-algebra $\mathbf{Pth}_{\boldsymbol{\mathcal{A}}^{(2)}}$; the ninth equality follows from the fact that $\mathrm{ip}^{(2,X)@}$ is a $\Sigma^{\boldsymbol{\mathcal{A}}^{(2)}}$-homomorphism, according to Definition~\ref{DDIp}; finally, the last equality recovers the value of the second-order Curry-Howard mapping at $\mathfrak{P}^{(2)}$, as we have discussed above.

The case~(1.1.2) of $\mathfrak{P}^{(2)}$ being a second-order path of length strictly greater than one containing a second-order echelon on its first step follows.

The case~(1.1), where $i=0$, follows.

For case~(1.2), i.e., if $i\neq 0$, since $\bb{\mathfrak{Q}^{(2)}\circ^{1\mathbf{Pth}_{\boldsymbol{\mathcal{A}}^{(2)}}}_{s}\mathfrak{P}^{(2)}}=\bb{\mathfrak{Q}^{(2)}}+\bb{\mathfrak{P}^{(2)}}$, then either (1.2.1) $i\in\bb{\mathfrak{P}^{(2)}}$ or (1.2.2) $i\in[\bb{\mathfrak{P}^{(2)}},\bb{\mathfrak{Q}^{(2)}\circ^{1\mathbf{Pth}_{\boldsymbol{\mathcal{A}}^{(2)}}}_{s}\mathfrak{P}^{(2)}}-1]$.

If~(1.2.1), i.e., if $i\neq 0$ and $i\in\bb{\mathfrak{P}^{(2)}}$, then $\mathfrak{P}^{(2)}$ is a second-order of length strictly greater than one containing a second-order echelon in a step different from the initial one. Then the value of the second-order Curry-Howard mapping at $\mathfrak{P}^{(2)}$ is given by
\[
\mathrm{CH}^{(2)}_{s}\left(
\mathfrak{P}^{(2)}
\right)
=
\mathrm{CH}^{(2)}_{s}\left(
\mathfrak{P}^{(2),i,\bb{\mathfrak{P}^{(2)}}-1}
\right)
\circ^{1\mathbf{T}_{\Sigma^{\boldsymbol{\mathcal{A}}^{(2)}}}(X)}_{s}
\mathrm{CH}^{(2)}_{s}\left(
\mathfrak{P}^{(2),0,i-1}
\right).
\]

Moreover, $(\mathfrak{Q}^{(2)}
\circ^{1\mathbf{Pth}_{\boldsymbol{\mathcal{A}}^{(2)}}}_{s}
\mathfrak{P}^{(2)})^{i,\bb{\mathfrak{Q}^{(2)}
\circ^{1\mathbf{Pth}_{\boldsymbol{\mathcal{A}}^{(2)}}}_{s}
\mathfrak{P}^{(2)}}-1}=\mathfrak{Q}^{(2)}\circ^{1\mathbf{Pth}_{\boldsymbol{\mathcal{A}}^{(2)}}}_{s}\mathfrak{P}^{(2),i,\bb{\mathfrak{P}^{(2))}-1}}$. Hence, the value of the second-order Curry-Howard mapping at $\mathfrak{Q}^{(2)}
\circ^{1\mathbf{Pth}_{\boldsymbol{\mathcal{A}}^{(2)}}}_{s}
\mathfrak{P}^{(2)}$ is given by 
\allowdisplaybreaks
\begin{multline*}
\mathrm{CH}^{(2)}_{s}\left(
\mathfrak{Q}^{(2)}
\circ^{1\mathbf{Pth}_{\boldsymbol{\mathcal{A}}^{(2)}}}_{s}
\mathfrak{P}^{(2)}
\right)
\\=
\mathrm{CH}^{(2)}_{s}\left(
\mathfrak{Q}^{(2)}\circ^{1\mathbf{Pth}_{\boldsymbol{\mathcal{A}}^{(2)}}}_{s}\mathfrak{P}^{(2),i,\bb{\mathfrak{P}^{(2))}-1}}
\right)
\circ^{1\mathbf{T}_{\Sigma^{\boldsymbol{\mathcal{A}}^{(2)}}}(X)}_{s}
\mathrm{CH}^{(2)}_{s}\left(
\mathfrak{P}^{(2),0,i}
\right).
\end{multline*}

Therefore, the following chain of equalities holds
\begin{flushleft}
$f_{s}\left(
\mathrm{ip}^{(2,X)@}_{s}\left(
\mathrm{CH}^{(2)}_{s}\left(
\mathfrak{Q}^{(2)}
\circ^{1\mathbf{Pth}_{\boldsymbol{\mathcal{A}}^{(2)}}}_{s}
\mathfrak{P}^{(2)}
\right)\right)
\right)$
\allowdisplaybreaks
\begin{align*}
&=f_{s}\left(
\mathrm{ip}^{(2,X)@}_{s}\left(
\mathrm{CH}^{(2)}_{s}\left(
\mathfrak{Q}^{(2)}\circ^{1\mathbf{Pth}_{\boldsymbol{\mathcal{A}}^{(2)}}}_{s}\mathfrak{P}^{(2),i,\bb{\mathfrak{P}^{(2))}-1}}
\right)
\right.\right.
\\&\qquad\qquad\qquad\qquad\qquad\qquad\qquad\qquad
\left.\left.
\circ^{1\mathbf{T}_{\Sigma^{\boldsymbol{\mathcal{A}}^{(2)}}}(X)}_{s}
\mathrm{CH}^{(2)}_{s}\left(
\mathfrak{P}^{(2),0,i-1}
\right)
\right)
\right)
\tag{1}
\\&=f_{s}\left(
\mathrm{ip}^{(2,X)@}_{s}\left(
\mathrm{CH}^{(2)}_{s}\left(
\mathfrak{Q}^{(2)}\circ^{1\mathbf{Pth}_{\boldsymbol{\mathcal{A}}^{(2)}}}_{s}\mathfrak{P}^{(2),i,\bb{\mathfrak{P}^{(2))}-1}}
\right)
\right)\right.
\\&\qquad\qquad\qquad\qquad\qquad\qquad
\left.
\circ^{1\mathbf{F}_{\Sigma^{\boldsymbol{\mathcal{A}}^{(2)}}}(\mathbf{Pth}_{\boldsymbol{\mathcal{A}}^{(2)}})}_{s}
\mathrm{ip}^{(2,X)@}_{s}\left(
\mathrm{CH}^{(2)}_{s}\left(
\mathfrak{P}^{(2),0,i-1}
\right)
\right)
\right)
\tag{2}
\\&=f_{s}\left(
\mathrm{ip}^{(2,X)@}_{s}\left(
\mathrm{CH}^{(2)}_{s}\left(
\mathfrak{Q}^{(2)}\circ^{1\mathbf{Pth}_{\boldsymbol{\mathcal{A}}^{(2)}}}_{s}\mathfrak{P}^{(2),i,\bb{\mathfrak{P}^{(2))}-1}}
\right)
\right)\right.
\\&\qquad\qquad\qquad\qquad\qquad\qquad\qquad\qquad
\left.
\circ^{1\mathbf{Pth}_{\boldsymbol{\mathcal{A}}^{(2)}}}_{s}
\mathrm{ip}^{(2,X)@}_{s}\left(
\mathrm{CH}^{(2)}_{s}\left(
\mathfrak{P}^{(2),0,i-1}
\right)
\right)
\right)
\tag{3}
\\&=f_{s}\left(
\mathrm{ip}^{(2,X)@}_{s}\left(
\mathrm{CH}^{(2)}_{s}\left(
\mathfrak{Q}^{(2)}\circ^{1\mathbf{Pth}_{\boldsymbol{\mathcal{A}}^{(2)}}}_{s}\mathfrak{P}^{(2),i,\bb{\mathfrak{P}^{(2))}-1}}
\right)
\right)\right)
\\&\qquad\qquad\qquad\qquad\qquad\qquad\qquad\qquad
\circ^{1\mathbf{B}}_{s}
f_{s}\left(
\mathrm{ip}^{(2,X)@}_{s}\left(
\mathrm{CH}^{(2)}_{s}\left(
\mathfrak{P}^{(2),0,i-1}
\right)
\right)
\right)
\tag{4}
\\&=\left(f_{s}\left(
\mathrm{ip}^{(2,X)@}_{s}\left(
\mathrm{CH}^{(2)}_{s}\left(
\mathfrak{Q}^{(2)}
\right)\right)\right)
\right.
\\&\qquad\qquad\qquad\quad
\left.
\circ^{1\mathbf{B}}_{s}
f_{s}\left(
\mathrm{ip}^{(2,X)@}_{s}\left(
\mathrm{CH}^{(2)}_{s}\left(
\mathfrak{P}^{(2),i,\bb{\mathfrak{P}^{(2))}-1}}
\right)
\right)\right)\right)
\\&\qquad\qquad\qquad\qquad\qquad\qquad\qquad\qquad
\circ^{1\mathbf{B}}_{s}
f_{s}\left(
\mathrm{ip}^{(2,X)@}_{s}\left(
\mathrm{CH}^{(2)}_{s}\left(
\mathfrak{P}^{(2),0,i-1}
\right)
\right)
\right)
\tag{5}
\\&=f_{s}\left(
\mathrm{ip}^{(2,X)@}_{s}\left(
\mathrm{CH}^{(2)}_{s}\left(
\mathfrak{Q}^{(2)}
\right)\right)\right)
\circ^{1\mathbf{B}}_{s}
\\&\qquad\qquad\qquad\quad
\left(
f_{s}\left(
\mathrm{ip}^{(2,X)@}_{s}\left(
\mathrm{CH}^{(2)}_{s}\left(
\mathfrak{P}^{(2),i,\bb{\mathfrak{P}^{(2))}-1}}
\right)
\right)\right)
\right.
\\&\qquad\qquad\qquad\qquad\qquad\qquad\qquad\quad
\circ^{1\mathbf{B}}_{s}
\left.
f_{s}\left(
\mathrm{ip}^{(2,X)@}_{s}\left(
\mathrm{CH}^{(2)}_{s}\left(
\mathfrak{P}^{(2),0,i-1}
\right)
\right)
\right)\right)
\tag{6}
\\&=f_{s}\left(
\mathrm{ip}^{(2,X)@}_{s}\left(
\mathrm{CH}^{(2)}_{s}\left(
\mathfrak{Q}^{(2)}
\right)\right)\right)
\circ^{1\mathbf{B}}_{s}
\\&\qquad\qquad\qquad\quad
\left(
f_{s}\left(
\mathrm{ip}^{(2,X)@}_{s}\left(
\mathrm{CH}^{(2)}_{s}\left(
\mathfrak{P}^{(2),i,\bb{\mathfrak{P}^{(2))}-1}}
\right)
\right)\right)
\right.
\\&\qquad\qquad\qquad\qquad\qquad\qquad\qquad
\circ^{1\mathbf{\mathbf{Pth}_{\boldsymbol{\mathcal{A}}^{(2)}}}}_{s}
\left.\left.
\mathrm{ip}^{(2,X)@}_{s}\left(
\mathrm{CH}^{(2)}_{s}\left(
\mathfrak{P}^{(2),0,i-1}
\right)\right)
\right)\right)
\tag{7}
\\&=f_{s}\left(
\mathrm{ip}^{(2,X)@}_{s}\left(
\mathrm{CH}^{(2)}_{s}\left(
\mathfrak{Q}^{(2)}
\right)\right)\right)
\circ^{1\mathbf{B}}_{s}
\\&\qquad\qquad\qquad\quad
\left(
f_{s}\left(
\mathrm{ip}^{(2,X)@}_{s}\left(
\mathrm{CH}^{(2)}_{s}\left(
\mathfrak{P}^{(2),i,\bb{\mathfrak{P}^{(2))}-1}}
\right)
\right)\right)
\right.
\\&\qquad\qquad\qquad\qquad\qquad
\circ^{1\mathbf{F}_{\Sigma^{\boldsymbol{\mathcal{A}}^{(2)}}}(\mathbf{\mathbf{Pth}_{\boldsymbol{\mathcal{A}}^{(2)}}})}_{s}
\left.\left.
\mathrm{ip}^{(2,X)@}_{s}\left(
\mathrm{CH}^{(2)}_{s}\left(
\mathfrak{P}^{(2),0,i-1}
\right)\right)
\right)\right)
\tag{8}
\\&=f_{s}\left(
\mathrm{ip}^{(2,X)@}_{s}\left(
\mathrm{CH}^{(2)}_{s}\left(
\mathfrak{Q}^{(2)}
\right)\right)\right)
\circ^{1\mathbf{B}}_{s}
\\&\qquad\qquad\qquad\qquad\qquad
\left(
f_{s}\left(
\mathrm{ip}^{(2,X)@}_{s}\left(
\mathrm{CH}^{(2)}_{s}\left(
\mathfrak{P}^{(2),i,\bb{\mathfrak{P}^{(2))}-1}}
\right)
\right.\right.\right.
\\&\qquad\qquad\qquad\qquad\qquad\qquad\qquad\qquad
\left.\left.\left.
\circ^{1\mathbf{T}_{\Sigma^{\boldsymbol{\mathcal{A}}^{(2)}}}(X)}_{s}
\mathrm{CH}^{(2)}_{s}\left(
\mathfrak{P}^{(2),0,i-1}
\right)
\right)
\right)
\right)
\tag{9}
\\&=
f_{s}\left(
\mathrm{ip}^{(2,X)@}_{s}\left(
\mathrm{CH}^{(2)}_{s}\left(
\mathfrak{Q}^{(2)}\right)\right)\right)
\circ^{1\mathbf{B}}_{s}
f_{s}\left(
\mathrm{ip}^{(2,X)@}_{s}\left(
\mathrm{CH}^{(2)}_{s}\left(
\mathfrak{P}^{(2)}
\right)\right)
\right).
\tag{10}
\end{align*}
\end{flushleft}

In the just stated chain of equalities, the first equality follows from the previous discussion on the value of the second-order Curry-Howard mapping at $\mathfrak{Q}^{(2)}\circ^{1\mathbf{Pth}_{\boldsymbol{\mathcal{A}}^{(2)}}}_{s}\mathfrak{P}^{(2)}$; the second equality follow from the fact that $\mathrm{ip}^{(2,X)@}$ is a $\Sigma^{\boldsymbol{\mathcal{A}}^{(2)}}$-homomorphism, according to Definition~\ref{DDIp}; the third equality follows from the fact that, by Proposition~\ref{PDIpDCH}, both $\mathrm{ip}^{(2,X)@}_{s}(
\mathrm{CH}^{(2)}_{s}(
\mathfrak{Q}^{(2)}\circ^{1\mathbf{Pth}_{\boldsymbol{\mathcal{A}}^{(2)}}}_{s}
\mathfrak{P}^{(2),i,\bb{\mathfrak{P}^{(2)}}-1}
))$ and $\mathrm{ip}^{(2,X)@}_{s}(
\mathrm{CH}^{(2)}_{s}(
\mathfrak{P}^{(2),0,i-1}
))$ are second-order paths. Thus, the interpretation of the operation symbol for the $1$-composition $\circ^{1}_{s}$ in the $\Sigma^{\boldsymbol{\mathcal{A}}^{(2)}}$-algebra $\mathbf{F}_{\Sigma^{\boldsymbol{\mathcal{A}}^{(2)}}}(\mathbf{Pth}_{\boldsymbol{\mathcal{A}}^{(2)}})$ becomes the corresponding interpretation of the $1$-composition $\circ^{1}_{s}$ in the many-sorted partial $\Sigma^{\boldsymbol{\mathcal{A}}^{(2)}}$-algebra $\mathbf{Pth}_{\boldsymbol{\mathcal{A}}^{(2)}}$; the fourth equality follows from the fact that, by assumption $f$ is a $\Sigma^{\boldsymbol{\mathcal{A}}^{(2)}}$-homomorphism; the fifth equality follows by induction. Let us note that the pair $(\mathfrak{Q}^{(2)}\circ^{1\mathbf{Pth}_{\boldsymbol{\mathcal{A}}^{(2)}}}_{s}
\mathfrak{P}^{(2),i,\bb{\mathfrak{P}^{(2)}}-1},s)$ $\prec_{\mathbf{Pth}_{\boldsymbol{\mathcal{A}}^{(2)}}}$-precedes $(\mathfrak{Q}^{(2)}\circ^{1\mathbf{Pth}_{\boldsymbol{\mathcal{A}}^{(2)}}}_{s}\mathfrak{P}^{(2)},s)$. Hence, we have that 
\begin{multline*}
f_{s}\left(
\mathrm{ip}^{(2,X)@}_{s}\left(
\mathrm{CH}^{(2)}_{s}\left(
\mathfrak{Q}^{(2)}
\circ^{1\mathbf{Pth}_{\boldsymbol{\mathcal{A}}^{(2)}}}_{s}
\mathfrak{P}^{(2),i,\bb{\mathfrak{P}^{(2)}}-1}
\right)\right)
\right)
\\=
f_{s}\left(
\mathrm{ip}^{(2,X)@}_{s}\left(
\mathrm{CH}^{(2)}_{s}\left(
\mathfrak{Q}^{(2)}
\right)\right)\right)
\circ^{1\mathbf{B}}_{s}
f_{s}\left(
\mathrm{ip}^{(2,X)@}_{s}\left(
\mathrm{CH}^{(2)}_{s}\left(
\mathfrak{P}^{(2),i,\bb{\mathfrak{P}^{(2)}}-1}
\right)\right)
\right);
\end{multline*}
the sixth equality follows from the fact that, since $\mathbf{B}$ is a partial many-sorted $\Sigma^{\boldsymbol{\mathcal{A}}^{(2)}}$-algebra in $\mathbf{PAlg}(\boldsymbol{\mathcal{E}}^{\boldsymbol{\mathcal{A}}^{(2)}})$, the $1$-composition is associative; the seventh equality follows from the fact that, by assumption, $f$ is a $\Sigma^{\boldsymbol{\mathcal{A}}^{(2)}}$-homomorphism; the eighth equality follows from the fact that,  by Proposition~\ref{PDIpDCH}, both $\mathrm{ip}^{(2,X)@}_{s}(
\mathrm{CH}^{(2)}_{s}(
\mathfrak{P}^{(2),i,\bb{\mathfrak{P}^{(2)}}-1}
))$ and $\mathrm{ip}^{(2,X)@}_{s}(
\mathrm{CH}^{(2)}_{s}(
\mathfrak{P}^{(2),0,i-1}
))$ are second-order paths. Thus, the interpretation of the operation symbol for the $1$-composition $\circ^{1}_{s}$ in the $\Sigma^{\boldsymbol{\mathcal{A}}^{(2)}}$-algebra $\mathbf{F}_{\Sigma^{\boldsymbol{\mathcal{A}}^{(2)}}}(\mathbf{Pth}_{\boldsymbol{\mathcal{A}}^{(2)}})$ becomes the corresponding interpretation of the $1$-composition $\circ^{1}_{s}$ in the many-sorted partial $\Sigma^{\boldsymbol{\mathcal{A}}^{(2)}}$-algebra $\mathbf{Pth}_{\boldsymbol{\mathcal{A}}^{(2)}}$; the ninth equality follows from the fact that $\mathrm{ip}^{(2,X)@}$ is a $\Sigma^{\boldsymbol{\mathcal{A}}^{(2)}}$-homomorphism, according to Definition~\ref{DDIp}; finally, the last equality recovers the value of the second-order Curry-Howard mapping at $\mathfrak{P}^{(2)}$, as we have discussed above.

The case~(1.2.1), where $i\in\bb{\mathfrak{P}^{(2)}}$ and $i\neq 0$ follows.

If~(1.2.2), i.e., if $i\neq 0$ and $i\in [\bb{\mathfrak{P}^{(2)}}, \bb{\mathfrak{Q}^{(2)}\circ^{1\mathbf{Pth}_{\boldsymbol{\mathcal{A}}^{(2)}}}_{s}\mathfrak{P}^{(2)}}-1]$ then $\mathfrak{Q}^{(2)}$ is a non-$(2,[1])$-identity second-order path containing a second-order echelon, whilst $\mathfrak{P}^{(2)}$ is an echelonless second-order path.

We will distinguish three cases according to whether  (1.2.2.1) $\mathfrak{Q}^{(2)}$ is a second-order echelon, (1.2.2.2) $\mathfrak{Q}^{(2)}$ is a second-order path of length strictly greater than one containing a second-order echelon on its first step or (1.2.2.3) $\mathfrak{Q}^{(2)}$ is a second-order path of length strictly greater than one containing a second-order echelon on a step different from zero. These cases can be proved following a similar argument to those three cases presented above. We leave the details for the interested reader.

If~(2), i.e., if $\mathfrak{Q}^{(2)}\circ^{1\mathbf{Pth}_{\boldsymbol{\mathcal{A}}^{(2)}}}_{s}\mathfrak{P}^{(2)}$ is an echelonless second-order path then, it could be the case that (2.1) $\mathfrak{Q}^{(2)}\circ^{1\mathbf{Pth}_{\boldsymbol{\mathcal{A}}^{(2)}}}_{s}\mathfrak{P}^{(2)}$ is an echelonless second-order path that is not head-constant, or (2.2) $\mathfrak{Q}^{(2)}\circ^{1\mathbf{Pth}_{\boldsymbol{\mathcal{A}}^{(2)}}}_{s}\mathfrak{P}^{(2)}$ is a head-constant echelonless second-order path that is not coherent, or (2.3) $\mathfrak{Q}^{(2)}\circ^{1\mathbf{Pth}_{\boldsymbol{\mathcal{A}}^{(2)}}}_{s}\mathfrak{P}^{(2)}$ is a coherent head-constant echelonless second-order path.

If~(2.1), let $i\in \bb{\mathfrak{Q}^{(2)}\circ^{1\mathbf{Pth}_{\boldsymbol{\mathcal{A}}^{(2)}}}_{s}\mathfrak{P}^{(2)}}$ be the greatest index for which $(\mathfrak{Q}^{(2)}\circ^{1\mathbf{Pth}_{\boldsymbol{\mathcal{A}}^{(2)}}}_{s}\mathfrak{P}^{(2)})^{0,i}$ is a head-constant second-order path. Since $\bb{\mathfrak{Q}^{(2)}\circ^{1\mathbf{Pth}_{\boldsymbol{\mathcal{A}}^{(2)}}}_{s}\mathfrak{P}^{(2)}}=\bb{\mathfrak{Q}^{(2)}}+\bb{\mathfrak{P}^{(2)}}$, we have that either (2.1.1) $i\in\bb{\mathfrak{P}^{(2)}}-1$, (2.1.2) $i=\bb{\mathfrak{P}^{(2)}}-1$, or (2.1.3) $i\in[\bb{\mathfrak{P}^{(2)}}, \bb{\mathfrak{Q}^{(2)}\circ^{1\mathbf{Pth}_{\boldsymbol{\mathcal{A}}^{(2)}}}_{s}\mathfrak{P}^{(2)}}-1]$.

If~(2.1.1), i.e., if we find ourselves in the case where $\mathfrak{Q}^{(2)}\circ^{1\mathbf{Pth}_{\boldsymbol{\mathcal{A}}^{(2)}}}_{s}\mathfrak{P}^{(2)}$ is an echelonless second-order path that is not head-constant and $i\in \bb{\mathfrak{P}^{(2)}}-1$ is the greatest index for which $(\mathfrak{Q}^{(2)}\circ^{1\mathbf{Pth}_{\boldsymbol{\mathcal{A}}^{(2)}}}_{s}\mathfrak{P}^{(2)})^{0,i}$ is a head-constant then, regarding the second-order paths $\mathfrak{Q}^{(2)}$ and $\mathfrak{P}^{(2)}$, we have that
\begin{itemize}
\item[(i)] $\mathfrak{P}^{(2)}$ is an echelonless second-orde path that is not head-constant and $i\in\bb{\mathfrak{P}^{(2)}}-1$ is the greatest index for which $\mathfrak{P}^{(2),0,i}$ is head-constant;
\item[(ii)] $\mathfrak{Q}^{(2)}$ is an echelonless second-order path.
\end{itemize}

From (i) and taking into account Definition~\ref{DDCH}, we have that the value of the second-order Curry-Howard mapping at $\mathfrak{P}^{(2)}$ is given by
\[
\mathrm{CH}^{(2)}_{s}\left(\mathfrak{P}^{(2)}\right)=
\mathrm{CH}^{(2)}_{s}\left(\mathfrak{P}^{(2),i+1,\bb{\mathfrak{P}^{(2)}}-1}\right)
\circ^{1\mathbf{T}_{\Sigma^{\boldsymbol{\mathcal{A}}^{(2)}}}(X)}_{s}
\mathrm{CH}^{(2)}_{s}\left(\mathfrak{P}^{(2),0,i}\right).
\]

Moreover, $(\mathfrak{Q}^{(2)}\circ^{1\mathbf{Pth}_{\boldsymbol{\mathcal{A}}^{(2)}}}_{s}\mathfrak{P}^{(2)})^{i+1,\bb{\mathfrak{Q}^{(2)}\circ^{1\mathbf{Pth}_{\boldsymbol{\mathcal{A}}^{(2)}}}_{s}\mathfrak{P}^{(2)}}-1}= \mathfrak{Q}^{(2)}\circ^{1\mathbf{Pth}_{\boldsymbol{\mathcal{A}}^{(2)}}}_{s}\mathfrak{P}^{(2),i+1,\bb{\mathfrak{P}^{(2)}}-1}$. Hence, the value of the second-order Curry-Howard mapping at  $\mathfrak{Q}^{(2)}
\circ^{1\mathbf{Pth}_{\boldsymbol{\mathcal{A}}^{(2)}}}_{s}
\mathfrak{P}^{(2)}$ is given by 
\allowdisplaybreaks
\begin{multline*}
\mathrm{CH}^{(2)}_{s}\left(
\mathfrak{Q}^{(2)}
\circ^{1\mathbf{Pth}_{\boldsymbol{\mathcal{A}}^{(2)}}}_{s}
\mathfrak{P}^{(2)}
\right)
\\=
\mathrm{CH}^{(2)}_{s}\left(
\mathfrak{Q}^{(2)}\circ^{1\mathbf{Pth}_{\boldsymbol{\mathcal{A}}^{(2)}}}_{s}\mathfrak{P}^{(2),i+1,\bb{\mathfrak{P}^{(2))}-1}}
\right)
\circ^{1\mathbf{T}_{\Sigma^{\boldsymbol{\mathcal{A}}^{(2)}}}(X)}_{s}
\mathrm{CH}^{(2)}_{s}\left(
\mathfrak{P}^{(2),0,i}
\right).
\end{multline*}

Therefore, the following chain of equalities holds
\begin{flushleft}
$f_{s}\left(
\mathrm{ip}^{(2,X)@}_{s}\left(
\mathrm{CH}^{(2)}_{s}\left(
\mathfrak{Q}^{(2)}
\circ^{1\mathbf{Pth}_{\boldsymbol{\mathcal{A}}^{(2)}}}_{s}
\mathfrak{P}^{(2)}
\right)\right)
\right)$
\allowdisplaybreaks
\begin{align*}
&=f_{s}\left(
\mathrm{ip}^{(2,X)@}_{s}\left(
\mathrm{CH}^{(2)}_{s}\left(
\mathfrak{Q}^{(2)}\circ^{1\mathbf{Pth}_{\boldsymbol{\mathcal{A}}^{(2)}}}_{s}\mathfrak{P}^{(2),i+1,\bb{\mathfrak{P}^{(2))}-1}}
\right)
\right.\right.
\\&\qquad\qquad\qquad\qquad\qquad\qquad\qquad\qquad
\left.\left.
\circ^{1\mathbf{T}_{\Sigma^{\boldsymbol{\mathcal{A}}^{(2)}}}(X)}_{s}
\mathrm{CH}^{(2)}_{s}\left(
\mathfrak{P}^{(2),0,i}
\right)
\right)
\right)
\tag{1}
\\&=f_{s}\left(
\mathrm{ip}^{(2,X)@}_{s}\left(
\mathrm{CH}^{(2)}_{s}\left(
\mathfrak{Q}^{(2)}\circ^{1\mathbf{Pth}_{\boldsymbol{\mathcal{A}}^{(2)}}}_{s}\mathfrak{P}^{(2),i+1,\bb{\mathfrak{P}^{(2))}-1}}
\right)
\right)\right.
\\&\qquad\qquad\qquad\qquad\qquad\qquad
\left.
\circ^{1\mathbf{F}_{\Sigma^{\boldsymbol{\mathcal{A}}^{(2)}}}(\mathbf{Pth}_{\boldsymbol{\mathcal{A}}^{(2)}})}_{s}
\mathrm{ip}^{(2,X)@}_{s}\left(
\mathrm{CH}^{(2)}_{s}\left(
\mathfrak{P}^{(2),0,i}
\right)
\right)
\right)
\tag{2}
\\&=f_{s}\left(
\mathrm{ip}^{(2,X)@}_{s}\left(
\mathrm{CH}^{(2)}_{s}\left(
\mathfrak{Q}^{(2)}\circ^{1\mathbf{Pth}_{\boldsymbol{\mathcal{A}}^{(2)}}}_{s}\mathfrak{P}^{(2),i+1,\bb{\mathfrak{P}^{(2))}-1}}
\right)
\right)\right.
\\&\qquad\qquad\qquad\qquad\qquad\qquad\qquad\qquad
\left.
\circ^{1\mathbf{Pth}_{\boldsymbol{\mathcal{A}}^{(2)}}}_{s}
\mathrm{ip}^{(2,X)@}_{s}\left(
\mathrm{CH}^{(2)}_{s}\left(
\mathfrak{P}^{(2),0,i}
\right)
\right)
\right)
\tag{3}
\\&=f_{s}\left(
\mathrm{ip}^{(2,X)@}_{s}\left(
\mathrm{CH}^{(2)}_{s}\left(
\mathfrak{Q}^{(2)}\circ^{1\mathbf{Pth}_{\boldsymbol{\mathcal{A}}^{(2)}}}_{s}\mathfrak{P}^{(2),i+1,\bb{\mathfrak{P}^{(2))}-1}}
\right)
\right)\right)
\\&\qquad\qquad\qquad\qquad\qquad\qquad\qquad\qquad
\circ^{1\mathbf{B}}_{s}
f_{s}\left(
\mathrm{ip}^{(2,X)@}_{s}\left(
\mathrm{CH}^{(2)}_{s}\left(
\mathfrak{P}^{(2),0,i}
\right)
\right)
\right)
\tag{4}
\\&=\left(f_{s}\left(
\mathrm{ip}^{(2,X)@}_{s}\left(
\mathrm{CH}^{(2)}_{s}\left(
\mathfrak{Q}^{(2)}
\right)\right)\right)
\right.
\\&\qquad\qquad\qquad\quad
\left.
\circ^{1\mathbf{B}}_{s}
f_{s}\left(
\mathrm{ip}^{(2,X)@}_{s}\left(
\mathrm{CH}^{(2)}_{s}\left(
\mathfrak{P}^{(2),i+1,\bb{\mathfrak{P}^{(2))}-1}}
\right)
\right)\right)\right)
\\&\qquad\qquad\qquad\qquad\qquad\qquad\qquad\qquad
\circ^{1\mathbf{B}}_{s}
f_{s}\left(
\mathrm{ip}^{(2,X)@}_{s}\left(
\mathrm{CH}^{(2)}_{s}\left(
\mathfrak{P}^{(2),0,i}
\right)
\right)
\right)
\tag{5}
\\&=f_{s}\left(
\mathrm{ip}^{(2,X)@}_{s}\left(
\mathrm{CH}^{(2)}_{s}\left(
\mathfrak{Q}^{(2)}
\right)\right)\right)
\circ^{1\mathbf{B}}_{s}
\\&\qquad\qquad\qquad\quad
\left(
f_{s}\left(
\mathrm{ip}^{(2,X)@}_{s}\left(
\mathrm{CH}^{(2)}_{s}\left(
\mathfrak{P}^{(2),i+1,\bb{\mathfrak{P}^{(2))}-1}}
\right)
\right)\right)
\right.
\\&\qquad\qquad\qquad\qquad\qquad\qquad\qquad\qquad
\circ^{1\mathbf{B}}_{s}
\left.
f_{s}\left(
\mathrm{ip}^{(2,X)@}_{s}\left(
\mathrm{CH}^{(2)}_{s}\left(
\mathfrak{P}^{(2),0,i}
\right)
\right)
\right)\right)
\tag{6}
\\&=f_{s}\left(
\mathrm{ip}^{(2,X)@}_{s}\left(
\mathrm{CH}^{(2)}_{s}\left(
\mathfrak{Q}^{(2)}
\right)\right)\right)
\circ^{1\mathbf{B}}_{s}
\\&\qquad\qquad\qquad\quad
\left(
f_{s}\left(
\mathrm{ip}^{(2,X)@}_{s}\left(
\mathrm{CH}^{(2)}_{s}\left(
\mathfrak{P}^{(2),i+1,\bb{\mathfrak{P}^{(2))}-1}}
\right)
\right)\right)
\right.
\\&\qquad\qquad\qquad\qquad\qquad\qquad\qquad
\circ^{1\mathbf{\mathbf{Pth}_{\boldsymbol{\mathcal{A}}^{(2)}}}}_{s}
\left.\left.
\mathrm{ip}^{(2,X)@}_{s}\left(
\mathrm{CH}^{(2)}_{s}\left(
\mathfrak{P}^{(2),0,i}
\right)\right)
\right)\right)
\tag{7}
\\&=f_{s}\left(
\mathrm{ip}^{(2,X)@}_{s}\left(
\mathrm{CH}^{(2)}_{s}\left(
\mathfrak{Q}^{(2)}
\right)\right)\right)
\circ^{1\mathbf{B}}_{s}
\\&\qquad\qquad\qquad\quad
\left(
f_{s}\left(
\mathrm{ip}^{(2,X)@}_{s}\left(
\mathrm{CH}^{(2)}_{s}\left(
\mathfrak{P}^{(2),i+1,\bb{\mathfrak{P}^{(2))}-1}}
\right)
\right)\right)
\right.
\\&\qquad\qquad\qquad\qquad\qquad\qquad
\circ^{1\mathbf{F}_{\Sigma^{\boldsymbol{\mathcal{A}}^{(2)}}}(\mathbf{\mathbf{Pth}_{\boldsymbol{\mathcal{A}}^{(2)}}})}_{s}
\left.\left.
\mathrm{ip}^{(2,X)@}_{s}\left(
\mathrm{CH}^{(2)}_{s}\left(
\mathfrak{P}^{(2),0,i}
\right)\right)
\right)\right)
\tag{8}
\\&=f_{s}\left(
\mathrm{ip}^{(2,X)@}_{s}\left(
\mathrm{CH}^{(2)}_{s}\left(
\mathfrak{Q}^{(2)}
\right)\right)\right)
\circ^{1\mathbf{B}}_{s}
\\&\qquad\qquad\qquad\qquad\qquad
\left(
f_{s}\left(
\mathrm{ip}^{(2,X)@}_{s}\left(
\mathrm{CH}^{(2)}_{s}\left(
\mathfrak{P}^{(2),i+1,\bb{\mathfrak{P}^{(2))}-1}}
\right)
\right.\right.\right.
\\&\qquad\qquad\qquad\qquad\qquad\qquad\qquad\qquad\qquad
\left.\left.\left.
\circ^{1\mathbf{T}_{\Sigma^{\boldsymbol{\mathcal{A}}^{(2)}}}(X)}_{s}
\mathrm{CH}^{(2)}_{s}\left(
\mathfrak{P}^{(2),0,i}
\right)
\right)
\right)
\right)
\tag{9}
\\&=
f_{s}\left(
\mathrm{ip}^{(2,X)@}_{s}\left(
\mathrm{CH}^{(2)}_{s}\left(
\mathfrak{Q}^{(2)}\right)\right)\right)
\circ^{1\mathbf{B}}_{s}
f_{s}\left(
\mathrm{ip}^{(2,X)@}_{s}\left(
\mathrm{CH}^{(2)}_{s}\left(
\mathfrak{P}^{(2)}
\right)\right)
\right).
\tag{10}
\end{align*}
\end{flushleft}

In the just stated chain of equalities, the first equality follows from the previous discussion on the value of the second-order Curry-Howard mapping at $\mathfrak{Q}^{(2)}\circ^{1\mathbf{Pth}_{\boldsymbol{\mathcal{A}}^{(2)}}}_{s}\mathfrak{P}^{(2)}$; the second equality follow from the fact that $\mathrm{ip}^{(2,X)@}$ is a $\Sigma^{\boldsymbol{\mathcal{A}}^{(2)}}$-homomorphism, according to Definition~\ref{DDIp}; the third equality follows from the fact that, by Proposition~\ref{PDIpDCH}, both $\mathrm{ip}^{(2,X)@}_{s}(
\mathrm{CH}^{(2)}_{s}(
\mathfrak{Q}^{(2)}\circ^{1\mathbf{Pth}_{\boldsymbol{\mathcal{A}}^{(2)}}}_{s}
\mathfrak{P}^{(2),i+1,\bb{\mathfrak{P}^{(2)}}-1}
))$ and $\mathrm{ip}^{(2,X)@}_{s}(
\mathrm{CH}^{(2)}_{s}(
\mathfrak{P}^{(2),0,i}
))$ are second-order paths. Thus, the interpretation of the operation symbol for the $1$-composition $\circ^{1}_{s}$ in the $\Sigma^{\boldsymbol{\mathcal{A}}^{(2)}}$-algebra $\mathbf{F}_{\Sigma^{\boldsymbol{\mathcal{A}}^{(2)}}}(\mathbf{Pth}_{\boldsymbol{\mathcal{A}}^{(2)}})$ becomes the corresponding interpretation of the $1$-composition $\circ^{1}_{s}$ in the many-sorted partial $\Sigma^{\boldsymbol{\mathcal{A}}^{(2)}}$-algebra $\mathbf{Pth}_{\boldsymbol{\mathcal{A}}^{(2)}}$; the fourth equality follows from the fact that, by assumption $f$ is a $\Sigma^{\boldsymbol{\mathcal{A}}^{(2)}}$-homomorphism; the fifth equality follows by induction. Let us note that the pair $(\mathfrak{Q}^{(2)}\circ^{1\mathbf{Pth}_{\boldsymbol{\mathcal{A}}^{(2)}}}_{s}
\mathfrak{P}^{(2),i+1,\bb{\mathfrak{P}^{(2)}}-1},s)$ $\prec_{\mathbf{Pth}_{\boldsymbol{\mathcal{A}}^{(2)}}}$-precedes $(\mathfrak{Q}^{(2)}\circ^{1\mathbf{Pth}_{\boldsymbol{\mathcal{A}}^{(2)}}}_{s}\mathfrak{P}^{(2)},s)$. Hence, we have that 
\begin{multline*}
f_{s}\left(
\mathrm{ip}^{(2,X)@}_{s}\left(
\mathrm{CH}^{(2)}_{s}\left(
\mathfrak{Q}^{(2)}
\circ^{1\mathbf{Pth}_{\boldsymbol{\mathcal{A}}^{(2)}}}_{s}
\mathfrak{P}^{(2),i+1,\bb{\mathfrak{P}^{(2)}}-1}
\right)\right)
\right)
\\=
f_{s}\left(
\mathrm{ip}^{(2,X)@}_{s}\left(
\mathrm{CH}^{(2)}_{s}\left(
\mathfrak{Q}^{(2)}
\right)\right)\right)
\circ^{1\mathbf{B}}_{s}
f_{s}\left(
\mathrm{ip}^{(2,X)@}_{s}\left(
\mathrm{CH}^{(2)}_{s}\left(
\mathfrak{P}^{(2),i+1,\bb{\mathfrak{P}^{(2)}}-1}
\right)\right)
\right);
\end{multline*}
the sixth equality follows from the fact that, since $\mathbf{B}$ is a partial many-sorted $\Sigma^{\boldsymbol{\mathcal{A}}^{(2)}}$-algebra in $\mathbf{PAlg}(\boldsymbol{\mathcal{E}}^{\boldsymbol{\mathcal{A}}^{(2)}})$, the $1$-composition is associative; the seventh equality follows from the fact that, by assumption, $f$ is a $\Sigma^{\boldsymbol{\mathcal{A}}^{(2)}}$-homomorphism; the eighth equality follows from the fact that,  by Proposition~\ref{PDIpDCH}, both $\mathrm{ip}^{(2,X)@}_{s}(
\mathrm{CH}^{(2)}_{s}(
\mathfrak{P}^{(2),i+1,\bb{\mathfrak{P}^{(2)}}-1}
))$ and $\mathrm{ip}^{(2,X)@}_{s}(
\mathrm{CH}^{(2)}_{s}(
\mathfrak{P}^{(2),0,i}
))$ are second-order paths. Thus, the interpretation of the operation symbol for the $1$-composition $\circ^{1}_{s}$ in the $\Sigma^{\boldsymbol{\mathcal{A}}^{(2)}}$-algebra $\mathbf{F}_{\Sigma^{\boldsymbol{\mathcal{A}}^{(2)}}}(\mathbf{Pth}_{\boldsymbol{\mathcal{A}}^{(2)}})$ becomes the corresponding interpretation of the $1$-composition $\circ^{1}_{s}$ in the many-sorted partial $\Sigma^{\boldsymbol{\mathcal{A}}^{(2)}}$-algebra $\mathbf{Pth}_{\boldsymbol{\mathcal{A}}^{(2)}}$; the ninth equality follows from the fact that $\mathrm{ip}^{(2,X)@}$ is a $\Sigma^{\boldsymbol{\mathcal{A}}^{(2)}}$-homomorphism, according to Definition~\ref{DDIp}; finally, the last equality recovers the value of the second-order Curry-Howard mapping at $\mathfrak{P}^{(2)}$, as we have discussed above.

The case~(2.1.1), where $i\in\bb{\mathfrak{P}^{(2)}}-1$ follows.

If~(2.1.2), i.e., if we find ourselves in the case where $\mathfrak{Q}^{(2)}\circ^{1\mathbf{Pth}_{\boldsymbol{\mathcal{A}}^{(2)}}}_{s}\mathfrak{P}^{(2)}$ is an echelonless second-order path that is not head constant and $i=\bb{\mathfrak{P}^{(2)}}-1$ is the greatest index for which $(\mathfrak{Q}^{(2)}\circ^{1\mathbf{Pth}_{\boldsymbol{\mathcal{A}}^{(2)}}}_{s}\mathfrak{P}^{(2)})^{0,i}$ is head-constant, then, regarding the second-order paths $\mathfrak{Q}^{(2)}$ and $\mathfrak{P}^{(2)}$, we have that
\begin{itemize}
\item[(i)] $\mathfrak{P}^{(2)}$ is a head-constant, echelonless second-order path.
\item[(ii)] $\mathfrak{Q}^{(2)}$ is an echelonless second-order path.
\end{itemize}

Moreover, $(\mathfrak{Q}^{(2)}\circ^{1\mathbf{Pth}_{\boldsymbol{\mathcal{A}}^{(2)}}}_{s}\mathfrak{P}^{(2)})^{i+1,\bb{\mathfrak{Q}^{(2)}\circ^{1\mathbf{Pth}_{\boldsymbol{\mathcal{A}}^{(2)}}}_{s}\mathfrak{P}^{(2)}}-1}= \mathfrak{Q}^{(2)}$. Hence, the value of the second-order Curry-Howard mapping at  $\mathfrak{Q}^{(2)}
\circ^{1\mathbf{Pth}_{\boldsymbol{\mathcal{A}}^{(2)}}}_{s}
\mathfrak{P}^{(2)}$ is given by 
\[
\mathrm{CH}^{(2)}_{s}\left(
\mathfrak{Q}^{(2)}
\circ^{1\mathbf{Pth}_{\boldsymbol{\mathcal{A}}^{(2)}}}_{s}
\mathfrak{P}^{(2)}
\right)
=
\mathrm{CH}^{(2)}_{s}\left(
\mathfrak{Q}^{(2)}
\right)
\circ^{1\mathbf{T}_{\Sigma^{\boldsymbol{\mathcal{A}}^{(2)}}}(X)}_{s}
\mathrm{CH}^{(2)}_{s}\left(
\mathfrak{P}^{(2)}
\right).
\]

Therefore, the following chain of equalities holds
\begin{flushleft}
$f_{s}\left(
\mathrm{ip}^{(2,X)@}_{s}\left(
\mathrm{CH}^{(2)}_{s}\left(
\mathfrak{Q}^{(2)}
\circ^{1\mathbf{Pth}_{\boldsymbol{\mathcal{A}}^{(2)}}}_{s}
\mathfrak{P}^{(2)}
\right)\right)
\right)$
\allowdisplaybreaks
\begin{align*}
&=f_{s}\left(
\mathrm{ip}^{(2,X)@}_{s}\left(
\mathrm{CH}^{(2)}_{s}\left(
\mathfrak{Q}^{(2)}
\right)
\circ^{1\mathbf{T}_{\Sigma^{\boldsymbol{\mathcal{A}}^{(2)}}}(X)}_{s}
\mathrm{CH}^{(2)}_{s}\left(
\mathfrak{P}^{(2)}
\right)
\right)
\right)
\tag{1}
\\&=f_{s}\left(
\mathrm{ip}^{(2,X)@}_{s}\left(
\mathrm{CH}^{(2)}_{s}\left(
\mathfrak{Q}^{(2)}
\right)
\right)
\circ^{1\mathbf{F}_{\Sigma^{\boldsymbol{\mathcal{A}}^{(2)}}}(\mathbf{Pth}_{\boldsymbol{\mathcal{A}}^{(2)}})}_{s}
\mathrm{ip}^{(2,X)@}_{s}\left(
\mathrm{CH}^{(2)}_{s}\left(
\mathfrak{P}^{(2)}
\right)
\right)
\right)
\tag{2}
\\&=f_{s}\left(
\mathrm{ip}^{(2,X)@}_{s}\left(
\mathrm{CH}^{(2)}_{s}\left(
\mathfrak{Q}^{(2)}
\right)
\right)
\circ^{1\mathbf{Pth}_{\boldsymbol{\mathcal{A}}^{(2)}}}_{s}
\mathrm{ip}^{(2,X)@}_{s}\left(
\mathrm{CH}^{(2)}_{s}\left(
\mathfrak{P}^{(2)}
\right)
\right)
\right)
\tag{3}
\\&=f_{s}\left(
\mathrm{ip}^{(2,X)@}_{s}\left(
\mathrm{CH}^{(2)}_{s}\left(
\mathfrak{Q}^{(2)}
\right)
\right)\right)
\circ^{1\mathbf{B}}_{s}
f_{s}\left(
\mathrm{ip}^{(2,X)@}_{s}\left(
\mathrm{CH}^{(2)}_{s}\left(
\mathfrak{P}^{(2)}
\right)
\right)\right).
\tag{4}
\end{align*}
\end{flushleft}

In the just stated chain of equalities, the first equality follows from the previous discussion on the value of the second-order Curry-Howard mapping at $\mathfrak{Q}^{(2)}\circ^{1\mathbf{Pth}_{\boldsymbol{\mathcal{A}}^{(2)}}}_{s}\mathfrak{P}^{(2)}$; the second equality follow from the fact that $\mathrm{ip}^{(2,X)@}$ is a $\Sigma^{\boldsymbol{\mathcal{A}}^{(2)}}$-homomorphism, according to Definition~\ref{DDIp}; the third equality follows from the fact that, by Proposition~\ref{PDIpDCH}, both $\mathrm{ip}^{(2,X)@}_{s}(
\mathrm{CH}^{(2)}_{s}(
\mathfrak{Q}^{(2)}\circ^{1\mathbf{Pth}_{\boldsymbol{\mathcal{A}}^{(2)}}}_{s}
\mathfrak{P}^{(2),i+1,\bb{\mathfrak{P}^{(2)}}-1}
))$ and $\mathrm{ip}^{(2,X)@}_{s}(
\mathrm{CH}^{(2)}_{s}(
\mathfrak{P}^{(2),0,i}
))$ are second-order paths. Thus, the interpretation of the operation symbol for the $1$-composition $\circ^{1}_{s}$ in the $\Sigma^{\boldsymbol{\mathcal{A}}^{(2)}}$-algebra $\mathbf{F}_{\Sigma^{\boldsymbol{\mathcal{A}}^{(2)}}}(\mathbf{Pth}_{\boldsymbol{\mathcal{A}}^{(2)}})$ becomes the corresponding interpretation of the $1$-composition $\circ^{1}_{s}$ in the many-sorted partial $\Sigma^{\boldsymbol{\mathcal{A}}^{(2)}}$-algebra $\mathbf{Pth}_{\boldsymbol{\mathcal{A}}^{(2)}}$; finally, the last equality follows from the fact that, by assumption $f$ is a $\Sigma^{\boldsymbol{\mathcal{A}}^{(2)}}$-homomorphism.

The case~(2.1.2), where $i=\bb{\mathfrak{P}^{(2)}}-1$ follows.

If~(2.1.3), i.e., if we find ourselves in the case where $\mathfrak{Q}^{(2)}\circ^{1\mathbf{Pth}_{\boldsymbol{\mathcal{A}}^{(2)}}}_{s}\mathfrak{P}^{(2)}$ is an echelonless second-order path that is not head-constant and $i\in [\bb{\mathfrak{P}^{(2)}}, \bb{\mathfrak{Q}^{(2)}\circ^{1\mathbf{Pth}_{\boldsymbol{\mathcal{A}}^{(2)}}}_{s}\mathfrak{P}^{(2)}}-1]$ is the greatest index for which $(\mathfrak{Q}^{(2)}\circ^{1\mathbf{Pth}_{\boldsymbol{\mathcal{A}}^{(2)}}}_{s}\mathfrak{P}^{(2)})^{0,i}$ is head-constant then, regarding the second-order paths $\mathfrak{Q}^{(2)}$ and $\mathfrak{P}^{(2)}$, we have that
\begin{itemize}
\item[(i)] $\mathfrak{P}^{(2)}$ is a head-constant echelonless second-order path;
\item[(ii)] $\mathfrak{Q}^{(2)}$ is an echelonless second-order path that is not head-constant and $i-\bb{\mathfrak{P}^{(2)}}\in \bb{\mathfrak{Q}^{(2)}}-1$ is the greatest index for which $\mathfrak{Q}^{(2),0,i-\bb{\mathfrak{P}^{(2)}}}$ is head-constant.
\end{itemize}

From (ii) and taking into account Definition~\ref{DDCH}, we have that the value of the second-order Curry-Howard mapping at $\mathfrak{Q}^{(2)}$ is given by
\[
\mathrm{CH}^{(2)}_{s}\left(\mathfrak{Q}^{(2)}\right)=
\mathrm{CH}^{(2)}_{s}\left(\mathfrak{Q}^{(2),i-\bb{\mathfrak{P}^{(2)}}+1,\bb{\mathfrak{Q}^{(2)}}-1}\right)
\circ^{1\mathbf{T}_{\Sigma^{\boldsymbol{\mathcal{A}}^{(2)}}}(X)}_{s}
\mathrm{CH}^{(2)}_{s}\left(\mathfrak{Q}^{(2),0,i-\bb{\mathfrak{P}^{(2)}}}\right).
\]

Moreover, $(\mathfrak{Q}^{(2)}\circ^{1\mathbf{Pth}_{\boldsymbol{\mathcal{A}}^{(2)}}}_{s}\mathfrak{P}^{(2)})^{i+1,\bb{\mathfrak{Q}^{(2)}\circ^{1\mathbf{Pth}_{\boldsymbol{\mathcal{A}}^{(2)}}}_{s}\mathfrak{P}^{(2)}}-1}= \mathfrak{Q}^{(2),0,i-\bb{\mathfrak{P}^{(2)}}}\circ^{1\mathbf{Pth}_{\boldsymbol{\mathcal{A}}^{(2)}}}_{s}\mathfrak{P}^{(2)}$. Hence, the value of the second-order Curry-Howard mapping at  $\mathfrak{Q}^{(2)}
\circ^{1\mathbf{Pth}_{\boldsymbol{\mathcal{A}}^{(2)}}}_{s}
\mathfrak{P}^{(2)}$ is given by 
\allowdisplaybreaks
\begin{multline*}
\mathrm{CH}^{(2)}_{s}\left(
\mathfrak{Q}^{(2)}
\circ^{1\mathbf{Pth}_{\boldsymbol{\mathcal{A}}^{(2)}}}_{s}
\mathfrak{P}^{(2)}
\right)
\\=
\mathrm{CH}^{(2)}_{s}\left(
\mathfrak{Q}^{(2),i-\bb{\mathfrak{P}^{(2)}}+1,\bb{\mathfrak{Q}^{(2)}}-1}
\right)
\circ^{1\mathbf{T}_{\Sigma^{\boldsymbol{\mathcal{A}}^{(2)}}}(X)}_{s}
\\
\mathrm{CH}^{(2)}_{s}\left(
\mathfrak{Q}^{(2),0,i-\bb{\mathfrak{P}^{(2)}}}\circ^{1\mathbf{Pth}_{\boldsymbol{\mathcal{A}}^{(2)}}}_{s}\mathfrak{P}^{(2)}
\right).
\end{multline*}

This case follows from a similar argument to that presented in Case~(2.1.1).

The case~(2.1.3), where $i\in [\bb{\mathfrak{P}^{(2)}}, \bb{\mathfrak{Q}^{(2)}\circ^{1\mathbf{Pth}_{\boldsymbol{\mathcal{A}}^{(2)}}}_{s}\mathfrak{P}^{(2)}}-1]$ follows.

This completes the case~(2.1).

If~(2.2), i.e., if $\mathfrak{Q}^{(2)}\circ^{1\mathbf{Pth}_{\boldsymbol{\mathcal{A}}^{(2)}}}_{s}\mathfrak{P}^{(2)}$ is a head-constant echelonless second-order path that is not coherent, then let $i\in \bb{\mathfrak{Q}^{(2)}\circ^{1\mathbf{Pth}_{\boldsymbol{\mathcal{A}}^{(2)}}}_{s}\mathfrak{P}^{(2)}}$ be the greatest index for which $(\mathfrak{Q}^{(2)}\circ^{1\mathbf{Pth}_{\boldsymbol{\mathcal{A}}^{(2)}}}_{s}\mathfrak{P}^{(2)})^{0,i}$ is a coherent head-constant second-order path. Since $\bb{\mathfrak{Q}^{(2)}\circ^{1\mathbf{Pth}_{\boldsymbol{\mathcal{A}}^{(2)}}}_{s}\mathfrak{P}^{(2)}}=\bb{\mathfrak{Q}^{(2)}}+\bb{\mathfrak{P}^{(2)}}$, we have that either (2.2.1) $i\in\bb{\mathfrak{P}^{(2)}}-1$, (2.2.2) $i=\bb{\mathfrak{P}^{(2)}}-1$, or (2.2.3) $i\in[\bb{\mathfrak{P}^{(2)}}, \bb{\mathfrak{Q}^{(2)}\circ^{1\mathbf{Pth}_{\boldsymbol{\mathcal{A}}^{(2)}}}_{s}\mathfrak{P}^{(2)}}-1]$.

If~(2.2.1), i.e., if we find ourselves in the case where $\mathfrak{Q}^{(2)}\circ^{1\mathbf{Pth}_{\boldsymbol{\mathcal{A}}^{(2)}}}_{s}\mathfrak{P}^{(2)}$ is a head-constant echelonless second-order path that is not coherent and $i\in \bb{\mathfrak{P}^{(2)}}-1$ is the greatest index for which $(\mathfrak{Q}^{(2)}\circ^{1\mathbf{Pth}_{\boldsymbol{\mathcal{A}}^{(2)}}}_{s}\mathfrak{P}^{(2)})^{0,i}$ is coherent head-constant then, regarding the second-order paths $\mathfrak{Q}^{(2)}$ and $\mathfrak{P}^{(2)}$, we have that
\begin{itemize}
\item[(i)] $\mathfrak{P}^{(2)}$ is a head-constant echelonless second-order path that is not coherent and $i\in\bb{\mathfrak{P}^{(2)}}-1$ is the greatest index for which $\mathfrak{P}^{(2),0,i}$ is coherent;
\item[(ii)] $\mathfrak{Q}^{(2)}$ is a head-constant echelonless second-order path.
\end{itemize}

From (i) and taking into account Definition~\ref{DDCH}, we have that the value of the second-order Curry-Howard mapping at $\mathfrak{P}^{(2)}$ is given by
\[
\mathrm{CH}^{(2)}_{s}\left(\mathfrak{P}^{(2)}\right)=
\mathrm{CH}^{(2)}_{s}\left(\mathfrak{P}^{(2),i+1,\bb{\mathfrak{P}^{(2)}}-1}\right)
\circ^{1\mathbf{T}_{\Sigma^{\boldsymbol{\mathcal{A}}^{(2)}}}(X)}_{s}
\mathrm{CH}^{(2)}_{s}\left(\mathfrak{P}^{(2),0,i}\right).
\]

Moreover, $(\mathfrak{Q}^{(2)}\circ^{1\mathbf{Pth}_{\boldsymbol{\mathcal{A}}^{(2)}}}_{s}\mathfrak{P}^{(2)})^{i+1,\bb{\mathfrak{Q}^{(2)}\circ^{1\mathbf{Pth}_{\boldsymbol{\mathcal{A}}^{(2)}}}_{s}\mathfrak{P}^{(2)}}-1}= \mathfrak{Q}^{(2)}\circ^{1\mathbf{Pth}_{\boldsymbol{\mathcal{A}}^{(2)}}}_{s}\mathfrak{P}^{(2),i+1,\bb{\mathfrak{P}^{(2)}}-1}$. Hence, the value of the second-order Curry-Howard mapping at  $\mathfrak{Q}^{(2)}
\circ^{1\mathbf{Pth}_{\boldsymbol{\mathcal{A}}^{(2)}}}_{s}
\mathfrak{P}^{(2)}$ is given by 
\allowdisplaybreaks
\begin{multline*}
\mathrm{CH}^{(2)}_{s}\left(
\mathfrak{Q}^{(2)}
\circ^{1\mathbf{Pth}_{\boldsymbol{\mathcal{A}}^{(2)}}}_{s}
\mathfrak{P}^{(2)}
\right)
\\=
\mathrm{CH}^{(2)}_{s}\left(
\mathfrak{Q}^{(2)}\circ^{1\mathbf{Pth}_{\boldsymbol{\mathcal{A}}^{(2)}}}_{s}\mathfrak{P}^{(2),i+1,\bb{\mathfrak{P}^{(2))}-1}}
\right)
\circ^{1\mathbf{T}_{\Sigma^{\boldsymbol{\mathcal{A}}^{(2)}}}(X)}_{s}
\mathrm{CH}^{(2)}_{s}\left(
\mathfrak{P}^{(2),0,i}
\right).
\end{multline*}

Therefore, the following chain of equalities holds
\begin{flushleft}
$f_{s}\left(
\mathrm{ip}^{(2,X)@}_{s}\left(
\mathrm{CH}^{(2)}_{s}\left(
\mathfrak{Q}^{(2)}
\circ^{1\mathbf{Pth}_{\boldsymbol{\mathcal{A}}^{(2)}}}_{s}
\mathfrak{P}^{(2)}
\right)\right)
\right)$
\allowdisplaybreaks
\begin{align*}
&=f_{s}\left(
\mathrm{ip}^{(2,X)@}_{s}\left(
\mathrm{CH}^{(2)}_{s}\left(
\mathfrak{Q}^{(2)}\circ^{1\mathbf{Pth}_{\boldsymbol{\mathcal{A}}^{(2)}}}_{s}\mathfrak{P}^{(2),i+1,\bb{\mathfrak{P}^{(2))}-1}}
\right)
\right.\right.
\\&\qquad\qquad\qquad\qquad\qquad\qquad\qquad\qquad
\left.\left.
\circ^{1\mathbf{T}_{\Sigma^{\boldsymbol{\mathcal{A}}^{(2)}}}(X)}_{s}
\mathrm{CH}^{(2)}_{s}\left(
\mathfrak{P}^{(2),0,i}
\right)
\right)
\right)
\tag{1}
\\&=f_{s}\left(
\mathrm{ip}^{(2,X)@}_{s}\left(
\mathrm{CH}^{(2)}_{s}\left(
\mathfrak{Q}^{(2)}\circ^{1\mathbf{Pth}_{\boldsymbol{\mathcal{A}}^{(2)}}}_{s}\mathfrak{P}^{(2),i+1,\bb{\mathfrak{P}^{(2))}-1}}
\right)
\right)\right.
\\&\qquad\qquad\qquad\qquad\qquad\qquad
\left.
\circ^{1\mathbf{F}_{\Sigma^{\boldsymbol{\mathcal{A}}^{(2)}}}(\mathbf{Pth}_{\boldsymbol{\mathcal{A}}^{(2)}})}_{s}
\mathrm{ip}^{(2,X)@}_{s}\left(
\mathrm{CH}^{(2)}_{s}\left(
\mathfrak{P}^{(2),0,i}
\right)
\right)
\right)
\tag{2}
\\&=f_{s}\left(
\mathrm{ip}^{(2,X)@}_{s}\left(
\mathrm{CH}^{(2)}_{s}\left(
\mathfrak{Q}^{(2)}\circ^{1\mathbf{Pth}_{\boldsymbol{\mathcal{A}}^{(2)}}}_{s}\mathfrak{P}^{(2),i+1,\bb{\mathfrak{P}^{(2))}-1}}
\right)
\right)\right.
\\&\qquad\qquad\qquad\qquad\qquad\qquad\qquad\qquad
\left.
\circ^{1\mathbf{Pth}_{\boldsymbol{\mathcal{A}}^{(2)}}}_{s}
\mathrm{ip}^{(2,X)@}_{s}\left(
\mathrm{CH}^{(2)}_{s}\left(
\mathfrak{P}^{(2),0,i}
\right)
\right)
\right)
\tag{3}
\\&=f_{s}\left(
\mathrm{ip}^{(2,X)@}_{s}\left(
\mathrm{CH}^{(2)}_{s}\left(
\mathfrak{Q}^{(2)}\circ^{1\mathbf{Pth}_{\boldsymbol{\mathcal{A}}^{(2)}}}_{s}\mathfrak{P}^{(2),i+1,\bb{\mathfrak{P}^{(2))}-1}}
\right)
\right)\right)
\\&\qquad\qquad\qquad\qquad\qquad\qquad\qquad\qquad
\circ^{1\mathbf{B}}_{s}
f_{s}\left(
\mathrm{ip}^{(2,X)@}_{s}\left(
\mathrm{CH}^{(2)}_{s}\left(
\mathfrak{P}^{(2),0,i}
\right)
\right)
\right)
\tag{4}
\\&=\left(f_{s}\left(
\mathrm{ip}^{(2,X)@}_{s}\left(
\mathrm{CH}^{(2)}_{s}\left(
\mathfrak{Q}^{(2)}
\right)\right)\right)
\right.
\\&\qquad\qquad\qquad\quad
\left.
\circ^{1\mathbf{B}}_{s}
f_{s}\left(
\mathrm{ip}^{(2,X)@}_{s}\left(
\mathrm{CH}^{(2)}_{s}\left(
\mathfrak{P}^{(2),i+1,\bb{\mathfrak{P}^{(2))}-1}}
\right)
\right)\right)\right)
\\&\qquad\qquad\qquad\qquad\qquad\qquad\qquad\qquad
\circ^{1\mathbf{B}}_{s}
f_{s}\left(
\mathrm{ip}^{(2,X)@}_{s}\left(
\mathrm{CH}^{(2)}_{s}\left(
\mathfrak{P}^{(2),0,i}
\right)
\right)
\right)
\tag{5}
\\&=f_{s}\left(
\mathrm{ip}^{(2,X)@}_{s}\left(
\mathrm{CH}^{(2)}_{s}\left(
\mathfrak{Q}^{(2)}
\right)\right)\right)
\circ^{1\mathbf{B}}_{s}
\\&\qquad\qquad\qquad\quad
\left(
f_{s}\left(
\mathrm{ip}^{(2,X)@}_{s}\left(
\mathrm{CH}^{(2)}_{s}\left(
\mathfrak{P}^{(2),i+1,\bb{\mathfrak{P}^{(2))}-1}}
\right)
\right)\right)
\right.
\\&\qquad\qquad\qquad\qquad\qquad\qquad\qquad\qquad
\circ^{1\mathbf{B}}_{s}
\left.
f_{s}\left(
\mathrm{ip}^{(2,X)@}_{s}\left(
\mathrm{CH}^{(2)}_{s}\left(
\mathfrak{P}^{(2),0,i}
\right)
\right)
\right)\right)
\tag{6}
\\&=f_{s}\left(
\mathrm{ip}^{(2,X)@}_{s}\left(
\mathrm{CH}^{(2)}_{s}\left(
\mathfrak{Q}^{(2)}
\right)\right)\right)
\circ^{1\mathbf{B}}_{s}
\\&\qquad\qquad\qquad\quad
\left(
f_{s}\left(
\mathrm{ip}^{(2,X)@}_{s}\left(
\mathrm{CH}^{(2)}_{s}\left(
\mathfrak{P}^{(2),i+1,\bb{\mathfrak{P}^{(2))}-1}}
\right)
\right)\right)
\right.
\\&\qquad\qquad\qquad\qquad\qquad\qquad\qquad
\circ^{1\mathbf{\mathbf{Pth}_{\boldsymbol{\mathcal{A}}^{(2)}}}}_{s}
\left.\left.
\mathrm{ip}^{(2,X)@}_{s}\left(
\mathrm{CH}^{(2)}_{s}\left(
\mathfrak{P}^{(2),0,i}
\right)\right)
\right)\right)
\tag{7}
\\&=f_{s}\left(
\mathrm{ip}^{(2,X)@}_{s}\left(
\mathrm{CH}^{(2)}_{s}\left(
\mathfrak{Q}^{(2)}
\right)\right)\right)
\circ^{1\mathbf{B}}_{s}
\\&\qquad\qquad\qquad\quad
\left(
f_{s}\left(
\mathrm{ip}^{(2,X)@}_{s}\left(
\mathrm{CH}^{(2)}_{s}\left(
\mathfrak{P}^{(2),i+1,\bb{\mathfrak{P}^{(2))}-1}}
\right)
\right)\right)
\right.
\\&\qquad\qquad\qquad\qquad\qquad\qquad
\circ^{1\mathbf{F}_{\Sigma^{\boldsymbol{\mathcal{A}}^{(2)}}}(\mathbf{\mathbf{Pth}_{\boldsymbol{\mathcal{A}}^{(2)}}})}_{s}
\left.\left.
\mathrm{ip}^{(2,X)@}_{s}\left(
\mathrm{CH}^{(2)}_{s}\left(
\mathfrak{P}^{(2),0,i}
\right)\right)
\right)\right)
\tag{8}
\\&=f_{s}\left(
\mathrm{ip}^{(2,X)@}_{s}\left(
\mathrm{CH}^{(2)}_{s}\left(
\mathfrak{Q}^{(2)}
\right)\right)\right)
\circ^{1\mathbf{B}}_{s}
\\&\qquad\qquad\qquad\qquad\qquad
\left(
f_{s}\left(
\mathrm{ip}^{(2,X)@}_{s}\left(
\mathrm{CH}^{(2)}_{s}\left(
\mathfrak{P}^{(2),i+1,\bb{\mathfrak{P}^{(2))}-1}}
\right)
\right.\right.\right.
\\&\qquad\qquad\qquad\qquad\qquad\qquad\qquad\qquad\qquad
\left.\left.\left.
\circ^{1\mathbf{T}_{\Sigma^{\boldsymbol{\mathcal{A}}^{(2)}}}(X)}_{s}
\mathrm{CH}^{(2)}_{s}\left(
\mathfrak{P}^{(2),0,i}
\right)
\right)
\right)
\right)
\tag{9}
\\&=
f_{s}\left(
\mathrm{ip}^{(2,X)@}_{s}\left(
\mathrm{CH}^{(2)}_{s}\left(
\mathfrak{Q}^{(2)}\right)\right)\right)
\circ^{1\mathbf{B}}_{s}
f_{s}\left(
\mathrm{ip}^{(2,X)@}_{s}\left(
\mathrm{CH}^{(2)}_{s}\left(
\mathfrak{P}^{(2)}
\right)\right)
\right).
\tag{10}
\end{align*}
\end{flushleft}

In the just stated chain of equalities, the first equality follows from the previous discussion on the value of the second-order Curry-Howard mapping at $\mathfrak{Q}^{(2)}\circ^{1\mathbf{Pth}_{\boldsymbol{\mathcal{A}}^{(2)}}}_{s}\mathfrak{P}^{(2)}$; the second equality follow from the fact that $\mathrm{ip}^{(2,X)@}$ is a $\Sigma^{\boldsymbol{\mathcal{A}}^{(2)}}$-homomorphism, according to Definition~\ref{DDIp}; the third equality follows from the fact that, by Proposition~\ref{PDIpDCH}, both $\mathrm{ip}^{(2,X)@}_{s}(
\mathrm{CH}^{(2)}_{s}(
\mathfrak{Q}^{(2)}\circ^{1\mathbf{Pth}_{\boldsymbol{\mathcal{A}}^{(2)}}}_{s}
\mathfrak{P}^{(2),i+1,\bb{\mathfrak{P}^{(2)}}-1}
))$ and $\mathrm{ip}^{(2,X)@}_{s}(
\mathrm{CH}^{(2)}_{s}(
\mathfrak{P}^{(2),0,i}
))$ are second-order paths. Thus, the interpretation of the operation symbol for the $1$-composition $\circ^{1}_{s}$ in the $\Sigma^{\boldsymbol{\mathcal{A}}^{(2)}}$-algebra $\mathbf{F}_{\Sigma^{\boldsymbol{\mathcal{A}}^{(2)}}}(\mathbf{Pth}_{\boldsymbol{\mathcal{A}}^{(2)}})$ becomes the corresponding interpretation of the $1$-composition $\circ^{1}_{s}$ in the many-sorted partial $\Sigma^{\boldsymbol{\mathcal{A}}^{(2)}}$-algebra $\mathbf{Pth}_{\boldsymbol{\mathcal{A}}^{(2)}}$; the fourth equality follows from the fact that, by assumption $f$ is a $\Sigma^{\boldsymbol{\mathcal{A}}^{(2)}}$-homomorphism; the fifth equality follows by induction. Let us note that the pair $(\mathfrak{Q}^{(2)}\circ^{1\mathbf{Pth}_{\boldsymbol{\mathcal{A}}^{(2)}}}_{s}
\mathfrak{P}^{(2),i+1,\bb{\mathfrak{P}^{(2)}}-1},s)$ $\prec_{\mathbf{Pth}_{\boldsymbol{\mathcal{A}}^{(2)}}}$-precedes $(\mathfrak{Q}^{(2)}\circ^{1\mathbf{Pth}_{\boldsymbol{\mathcal{A}}^{(2)}}}_{s}\mathfrak{P}^{(2)},s)$. Hence, we have that 
\begin{multline*}
f_{s}\left(
\mathrm{ip}^{(2,X)@}_{s}\left(
\mathrm{CH}^{(2)}_{s}\left(
\mathfrak{Q}^{(2)}
\circ^{1\mathbf{Pth}_{\boldsymbol{\mathcal{A}}^{(2)}}}_{s}
\mathfrak{P}^{(2),i+1,\bb{\mathfrak{P}^{(2)}}-1}
\right)\right)
\right)
\\=
f_{s}\left(
\mathrm{ip}^{(2,X)@}_{s}\left(
\mathrm{CH}^{(2)}_{s}\left(
\mathfrak{Q}^{(2)}
\right)\right)\right)
\circ^{1\mathbf{B}}_{s}
f_{s}\left(
\mathrm{ip}^{(2,X)@}_{s}\left(
\mathrm{CH}^{(2)}_{s}\left(
\mathfrak{P}^{(2),i+1,\bb{\mathfrak{P}^{(2)}}-1}
\right)\right)
\right);
\end{multline*}
the sixth equality follows from the fact that, since $\mathbf{B}$ is a partial many-sorted $\Sigma^{\boldsymbol{\mathcal{A}}^{(2)}}$-algebra in $\mathbf{PAlg}(\boldsymbol{\mathcal{E}}^{\boldsymbol{\mathcal{A}}^{(2)}})$, the $1$-composition is associative; the seventh equality follows from the fact that, by assumption, $f$ is a $\Sigma^{\boldsymbol{\mathcal{A}}^{(2)}}$-homomorphism; the eighth equality follows from the fact that,  by Proposition~\ref{PDIpDCH}, both $\mathrm{ip}^{(2,X)@}_{s}(
\mathrm{CH}^{(2)}_{s}(
\mathfrak{P}^{(2),i+1,\bb{\mathfrak{P}^{(2)}}-1}
))$ and $\mathrm{ip}^{(2,X)@}_{s}(
\mathrm{CH}^{(2)}_{s}(
\mathfrak{P}^{(2),0,i}
))$ are second-order paths. Thus, the interpretation of the operation symbol for the $1$-composition $\circ^{1}_{s}$ in the $\Sigma^{\boldsymbol{\mathcal{A}}^{(2)}}$-algebra $\mathbf{F}_{\Sigma^{\boldsymbol{\mathcal{A}}^{(2)}}}(\mathbf{Pth}_{\boldsymbol{\mathcal{A}}^{(2)}})$ becomes the corresponding interpretation of the $1$-composition $\circ^{1}_{s}$ in the many-sorted partial $\Sigma^{\boldsymbol{\mathcal{A}}^{(2)}}$-algebra $\mathbf{Pth}_{\boldsymbol{\mathcal{A}}^{(2)}}$; the ninth equality follows from the fact that $\mathrm{ip}^{(2,X)@}$ is a $\Sigma^{\boldsymbol{\mathcal{A}}^{(2)}}$-homomorphism, according to Definition~\ref{DDIp}; finally, the last equality recovers the value of the second-order Curry-Howard mapping at $\mathfrak{P}^{(2)}$, as we have discussed above.

The case~(2.2.1), where $i\in\bb{\mathfrak{P}^{(2)}}-1$ follows.

If~(2.2.2), i.e., if we find ourselves in the case where $\mathfrak{Q}^{(2)}\circ^{1\mathbf{Pth}_{\boldsymbol{\mathcal{A}}^{(2)}}}_{s}\mathfrak{P}^{(2)}$ is a head-constant echelonless second-order path that is not coherent and $i=\bb{\mathfrak{P}^{(2)}}-1$ is the greatest index for which $(\mathfrak{Q}^{(2)}\circ^{1\mathbf{Pth}_{\boldsymbol{\mathcal{A}}^{(2)}}}_{s}\mathfrak{P}^{(2)})^{0,i}$ is coherent, then, regarding the second-order paths $\mathfrak{Q}^{(2)}$ and $\mathfrak{P}^{(2)}$, we have that
\begin{itemize}
\item[(i)] $\mathfrak{P}^{(2)}$ is a coherent head-constant echelonless second-order path.
\item[(ii)] $\mathfrak{Q}^{(2)}$ is a head-constant echelonless second-order path.
\end{itemize}

Moreover, $(\mathfrak{Q}^{(2)}\circ^{1\mathbf{Pth}_{\boldsymbol{\mathcal{A}}^{(2)}}}_{s}\mathfrak{P}^{(2)})^{i+1,\bb{\mathfrak{Q}^{(2)}\circ^{1\mathbf{Pth}_{\boldsymbol{\mathcal{A}}^{(2)}}}_{s}\mathfrak{P}^{(2)}}-1}= \mathfrak{Q}^{(2)}$. Hence, the value of the second-order Curry-Howard mapping at  $\mathfrak{Q}^{(2)}
\circ^{1\mathbf{Pth}_{\boldsymbol{\mathcal{A}}^{(2)}}}_{s}
\mathfrak{P}^{(2)}$ is given by 
\[
\mathrm{CH}^{(2)}_{s}\left(
\mathfrak{Q}^{(2)}
\circ^{1\mathbf{Pth}_{\boldsymbol{\mathcal{A}}^{(2)}}}_{s}
\mathfrak{P}^{(2)}
\right)
=
\mathrm{CH}^{(2)}_{s}\left(
\mathfrak{Q}^{(2)}
\right)
\circ^{1\mathbf{T}_{\Sigma^{\boldsymbol{\mathcal{A}}^{(2)}}}(X)}_{s}
\mathrm{CH}^{(2)}_{s}\left(
\mathfrak{P}^{(2)}
\right).
\]

Therefore, the following chain of equalities holds
\begin{flushleft}
$f_{s}\left(
\mathrm{ip}^{(2,X)@}_{s}\left(
\mathrm{CH}^{(2)}_{s}\left(
\mathfrak{Q}^{(2)}
\circ^{1\mathbf{Pth}_{\boldsymbol{\mathcal{A}}^{(2)}}}_{s}
\mathfrak{P}^{(2)}
\right)\right)
\right)$
\allowdisplaybreaks
\begin{align*}
&=f_{s}\left(
\mathrm{ip}^{(2,X)@}_{s}\left(
\mathrm{CH}^{(2)}_{s}\left(
\mathfrak{Q}^{(2)}
\right)
\circ^{1\mathbf{T}_{\Sigma^{\boldsymbol{\mathcal{A}}^{(2)}}}(X)}_{s}
\mathrm{CH}^{(2)}_{s}\left(
\mathfrak{P}^{(2)}
\right)
\right)
\right)
\tag{1}
\\&=f_{s}\left(
\mathrm{ip}^{(2,X)@}_{s}\left(
\mathrm{CH}^{(2)}_{s}\left(
\mathfrak{Q}^{(2)}
\right)
\right)
\circ^{1\mathbf{F}_{\Sigma^{\boldsymbol{\mathcal{A}}^{(2)}}}(\mathbf{Pth}_{\boldsymbol{\mathcal{A}}^{(2)}})}_{s}
\mathrm{ip}^{(2,X)@}_{s}\left(
\mathrm{CH}^{(2)}_{s}\left(
\mathfrak{P}^{(2)}
\right)
\right)
\right)
\tag{2}
\\&=f_{s}\left(
\mathrm{ip}^{(2,X)@}_{s}\left(
\mathrm{CH}^{(2)}_{s}\left(
\mathfrak{Q}^{(2)}
\right)
\right)
\circ^{1\mathbf{Pth}_{\boldsymbol{\mathcal{A}}^{(2)}}}_{s}
\mathrm{ip}^{(2,X)@}_{s}\left(
\mathrm{CH}^{(2)}_{s}\left(
\mathfrak{P}^{(2)}
\right)
\right)
\right)
\tag{3}
\\&=f_{s}\left(
\mathrm{ip}^{(2,X)@}_{s}\left(
\mathrm{CH}^{(2)}_{s}\left(
\mathfrak{Q}^{(2)}
\right)
\right)\right)
\circ^{1\mathbf{B}}_{s}
f_{s}\left(
\mathrm{ip}^{(2,X)@}_{s}\left(
\mathrm{CH}^{(2)}_{s}\left(
\mathfrak{P}^{(2)}
\right)
\right)\right).
\tag{4}
\end{align*}
\end{flushleft}

In the just stated chain of equalities, the first equality follows from the previous discussion on the value of the second-order Curry-Howard mapping at $\mathfrak{Q}^{(2)}\circ^{1\mathbf{Pth}_{\boldsymbol{\mathcal{A}}^{(2)}}}_{s}\mathfrak{P}^{(2)}$; the second equality follow from the fact that $\mathrm{ip}^{(2,X)@}$ is a $\Sigma^{\boldsymbol{\mathcal{A}}^{(2)}}$-homomorphism, according to Definition~\ref{DDIp}; the third equality follows from the fact that, by Proposition~\ref{PDIpDCH}, both $\mathrm{ip}^{(2,X)@}_{s}(
\mathrm{CH}^{(2)}_{s}(
\mathfrak{Q}^{(2)}\circ^{1\mathbf{Pth}_{\boldsymbol{\mathcal{A}}^{(2)}}}_{s}
\mathfrak{P}^{(2),i+1,\bb{\mathfrak{P}^{(2)}}-1}
))$ and $\mathrm{ip}^{(2,X)@}_{s}(
\mathrm{CH}^{(2)}_{s}(
\mathfrak{P}^{(2),0,i}
))$ are second-order paths. Thus, the interpretation of the operation symbol for the $1$-composition $\circ^{1}_{s}$ in the $\Sigma^{\boldsymbol{\mathcal{A}}^{(2)}}$-algebra $\mathbf{F}_{\Sigma^{\boldsymbol{\mathcal{A}}^{(2)}}}(\mathbf{Pth}_{\boldsymbol{\mathcal{A}}^{(2)}})$ becomes the corresponding interpretation of the $1$-composition $\circ^{1}_{s}$ in the many-sorted partial $\Sigma^{\boldsymbol{\mathcal{A}}^{(2)}}$-algebra $\mathbf{Pth}_{\boldsymbol{\mathcal{A}}^{(2)}}$; finally, the last equality follows from the fact that, by assumption $f$ is a $\Sigma^{\boldsymbol{\mathcal{A}}^{(2)}}$-homomorphism.

The case~(2.2.2), where $i=\bb{\mathfrak{P}^{(2)}}-1$ follows.

If~(2.2.3), i.e., if we find ourselves in the case where $\mathfrak{Q}^{(2)}\circ^{1\mathbf{Pth}_{\boldsymbol{\mathcal{A}}^{(2)}}}_{s}\mathfrak{P}^{(2)}$ is a head-constant echelonless second-order path that is not coherent and $i\in [\bb{\mathfrak{P}^{(2)}}, \bb{\mathfrak{Q}^{(2)}\circ^{1\mathbf{Pth}_{\boldsymbol{\mathcal{A}}^{(2)}}}_{s}\mathfrak{P}^{(2)}}-1]$ is the greatest index for which $(\mathfrak{Q}^{(2)}\circ^{1\mathbf{Pth}_{\boldsymbol{\mathcal{A}}^{(2)}}}_{s}\mathfrak{P}^{(2)})^{0,i}$ is coherent then, regarding the second-order paths $\mathfrak{Q}^{(2)}$ and $\mathfrak{P}^{(2)}$, we have that
\begin{itemize}
\item[(i)] $\mathfrak{P}^{(2)}$ is a coherent head-constant echelonless second-order path;
\item[(ii)] $\mathfrak{Q}^{(2)}$ is a head-constant echelonless second-order path that is not coherent and $i-\bb{\mathfrak{P}^{(2)}}\in \bb{\mathfrak{Q}^{(2)}}-1$ is the greatest index for which $\mathfrak{Q}^{(2),0,i-\bb{\mathfrak{P}^{(2)}}}$ is coherent.
\end{itemize}

From (ii) and taking into account Definition~\ref{DDCH}, we have that the value of the second-order Curry-Howard mapping at $\mathfrak{Q}^{(2)}$ is given by
\[
\mathrm{CH}^{(2)}_{s}\left(\mathfrak{Q}^{(2)}\right)=
\mathrm{CH}^{(2)}_{s}\left(\mathfrak{Q}^{(2),i-\bb{\mathfrak{P}^{(2)}}+1,\bb{\mathfrak{Q}^{(2)}}-1}\right)
\circ^{1\mathbf{T}_{\Sigma^{\boldsymbol{\mathcal{A}}^{(2)}}}(X)}_{s}
\mathrm{CH}^{(2)}_{s}\left(\mathfrak{Q}^{(2),0,i-\bb{\mathfrak{P}^{(2)}}}\right).
\]

Moreover, $(\mathfrak{Q}^{(2)}\circ^{1\mathbf{Pth}_{\boldsymbol{\mathcal{A}}^{(2)}}}_{s}\mathfrak{P}^{(2)})^{i+1,\bb{\mathfrak{Q}^{(2)}\circ^{1\mathbf{Pth}_{\boldsymbol{\mathcal{A}}^{(2)}}}_{s}\mathfrak{P}^{(2)}}-1}= \mathfrak{Q}^{(2),0,i-\bb{\mathfrak{P}^{(2)}}}\circ^{1\mathbf{Pth}_{\boldsymbol{\mathcal{A}}^{(2)}}}_{s}\mathfrak{P}^{(2)}$. Hence, the value of the second-order Curry-Howard mapping at  $\mathfrak{Q}^{(2)}
\circ^{1\mathbf{Pth}_{\boldsymbol{\mathcal{A}}^{(2)}}}_{s}
\mathfrak{P}^{(2)}$ is given by 
\allowdisplaybreaks
\begin{multline*}
\mathrm{CH}^{(2)}_{s}\left(
\mathfrak{Q}^{(2)}
\circ^{1\mathbf{Pth}_{\boldsymbol{\mathcal{A}}^{(2)}}}_{s}
\mathfrak{P}^{(2)}
\right)
\\=
\mathrm{CH}^{(2)}_{s}\left(
\mathfrak{Q}^{(2),i-\bb{\mathfrak{P}^{(2)}}+1,\bb{\mathfrak{Q}^{(2)}}-1}
\right)
\circ^{1\mathbf{T}_{\Sigma^{\boldsymbol{\mathcal{A}}^{(2)}}}(X)}_{s}
\\
\mathrm{CH}^{(2)}_{s}\left(
\mathfrak{Q}^{(2),0,i-\bb{\mathfrak{P}^{(2)}}}\circ^{1\mathbf{Pth}_{\boldsymbol{\mathcal{A}}^{(2)}}}_{s}\mathfrak{P}^{(2)}
\right).
\end{multline*}

This case follows from a similar argument to that presented in Case~(2.2.1).

The case~(2.2.3), where $i\in [\bb{\mathfrak{P}^{(2)}}, \bb{\mathfrak{Q}^{(2)}\circ^{1\mathbf{Pth}_{\boldsymbol{\mathcal{A}}^{(2)}}}_{s}\mathfrak{P}^{(2)}}-1]$ follows.

This completes the case~(2.2).

If~(2.3), i.e., if we find ourselves in the case  where  $\mathfrak{Q}^{(2)}\circ^{1\mathbf{Pth}_{\boldsymbol{\mathcal{A}}^{(2)}}}_{s}\mathfrak{P}^{(2)}$ is a coherent head-constant echelonless second-order path then, regarding the second-order paths $\mathfrak{Q}^{(2)}$ and $\mathfrak{P}^{(2)}$, we have that
\begin{itemize}
\item[(i)] $\mathfrak{P}^{(2)}$ is a coherent head-constant echelonless second-order path;
\item[(ii)] $\mathfrak{Q}^{(2)}$ is a coherent head-constant echelonless second-order path;
\end{itemize}

Therefore, for a unique word $\mathbf{s}\in S^{\star}-\{\lambda\}$ and a unique operation symbol $\tau\in\Sigma^{\boldsymbol{\mathcal{A}}}_{\mathbf{s},s}$, the family of first-order translations occurring in $\mathfrak{P}^{(2)}$ is a family of first-order translations of type $\tau$.

Let $(\mathfrak{P}^{(2)}_{j})_{j\in\bb{\mathbf{s}}}$ be the family of second-order paths we can extract from $\mathfrak{P}^{(2)}$ in virtue of Lemma~\ref{LDPthExtract}. Then according to Definition~\ref{DDCH}, we have that the value of the second-order Curry-Howard mapping at $\mathfrak{P}^{(2)}$ is given by
\[
\mathrm{CH}^{(2)}_{s}\left(
\mathfrak{P}^{(2)}
\right)
=
\tau^{\mathbf{T}_{\Sigma^{\boldsymbol{\mathcal{A}}^{(2)}}}(X)}
\left(\left(
\mathrm{CH}^{(2)}_{s_{j}}\left(
\mathfrak{P}^{(2)}_{j}
\right)\right)_{j\in\bb{\mathbf{s}}}
\right).
\]

Since (ii), we have that, for the unique word $\mathbf{s}\in S^{\star}-\{\lambda\}$ and the unique operation symbol $\tau\in\Sigma^{\boldsymbol{\mathcal{A}}}_{\mathbf{s},s}$,the family of first-order translations occurring in $\mathfrak{Q}^{(2)}$ is a family of first-order translations of type $\tau$.

Note that the operation symbol $\tau$ is the same as in case (i), since $\mathfrak{Q}^{(2)}\circ^{1\mathbf{Pth}_{\boldsymbol{\mathcal{A}}^{(2)}}}_{s}\mathfrak{P}^{(2)}$ is head-constant by hypothesis.

Let $(\mathfrak{Q}^{(2)}_{j})_{j\in\bb{\mathbf{s}}}$ be the family of second-order paths we can extract from $\mathfrak{Q}^{(2)}$ in virtue of Lemma~\ref{LDPthExtract}. Then according to Definition~\ref{DDCH}, we have that the value of the second-order Curry-Howard mapping at $\mathfrak{Q}^{(2)}$ is given by
\[
\mathrm{CH}^{(2)}_{s}\left(
\mathfrak{Q}^{(2)}
\right)
=
\tau^{\mathbf{T}_{\Sigma^{\boldsymbol{\mathcal{A}}^{(2)}}}(X)}
\left(\left(
\mathrm{CH}^{(2)}_{s_{j}}\left(
\mathfrak{Q}^{(2)}_{j}
\right)\right)_{j\in\bb{\mathbf{s}}}
\right).
\]

Let $((\mathfrak{Q}^{(2)}\circ^{1\mathbf{Pth}_{\boldsymbol{\mathcal{A}}^{(2)}}}_{s}\mathfrak{P}^{(2)})_{j})_{j\in\bb{\mathbf{s}}}$ be the family of second-order paths in $\mathrm{Pth}_{\boldsymbol{\mathcal{A}}^{(2)},\mathbf{s}}$ in virtue of Lemma~\ref{LDPthExtract}. Then, for every $j\in\bb{\mathbf{s}}$, we have that
\[
\left(
\mathfrak{Q}^{(2)}\circ^{1\mathbf{Pth}_{\boldsymbol{\mathcal{A}}^{(2)}}}_{s}\mathfrak{P}^{(2)}
\right)_{j}
=
\mathfrak{Q}^{(2)}_{j}
\circ^{1\mathbf{Pth}_{\boldsymbol{\mathcal{A}}^{(2)}}}_{s_{j}}
\mathfrak{P}^{(2)}_{j}.
\]

Thus, the value of the second-order Curry-Howard mapping at $\mathfrak{Q}^{(2)}\circ^{1\mathbf{Pth}_{\boldsymbol{\mathcal{A}}^{(2)}}}_{s}\mathfrak{P}^{(2)}$ is given by
\[
\mathrm{CH}^{(2)}_{s}\left(
\mathfrak{Q}^{(2)}\circ^{1\mathbf{Pth}_{\boldsymbol{\mathcal{A}}^{(2)}}}_{s}\mathfrak{P}^{(2)}
\right)
=
\tau^{\mathbf{T}_{\Sigma^{\boldsymbol{\mathcal{A}}^{(2)}}}(X)}
\left(\left(
\mathrm{CH}^{(2)}_{s_{j}}\left(
\mathfrak{Q}^{(2)}_{j}
\circ^{1\mathbf{Pth}_{\boldsymbol{\mathcal{A}}^{(2)}}}_{s_{j}}
\mathfrak{P}^{(2)}_{j}
\right)
\right)_{j\in\bb{\mathbf{s}}}
\right).
\]

Therefore, the following chain of equalities holds
\begin{flushleft}
$f_{s}\left(
\mathrm{ip}^{(2,X)@}_{s}\left(
\mathrm{CH}^{(2)}_{s}\left(
\mathfrak{Q}^{(2)}
\circ^{1\mathbf{Pth}_{\boldsymbol{\mathcal{A}}^{(2)}}}_{s}
\mathfrak{P}^{(2)}
\right)\right)
\right)$
\allowdisplaybreaks
\begin{align*}
&=
f_{s}\left(
\mathrm{ip}^{(2,X)@}_{s}\left(
\tau^{\mathbf{T}_{\Sigma^{\boldsymbol{\mathcal{A}}^{(2)}}}(X)}
\left(\left(
\mathrm{CH}^{(2)}_{s_{j}}\left(
\mathfrak{Q}^{(2)}_{j}
\circ^{1\mathbf{Pth}_{\boldsymbol{\mathcal{A}}^{(2)}}}_{s_{j}}
\mathfrak{P}^{(2)}_{j}
\right)
\right)_{j\in\bb{\mathbf{s}}}
\right)
\right)\right)
\tag{1}
\\&=
f_{s}\left(
\tau^{\mathbf{F}_{\Sigma^{\boldsymbol{\mathcal{A}}^{(2)}}}(\mathbf{Pth}_{\boldsymbol{\mathcal{A}}^{(2)}})}
\left(\left(
\mathrm{ip}^{(2,X)@}_{s_{j}}\left(
\mathrm{CH}^{(2)}_{s_{j}}\left(
\mathfrak{Q}^{(2)}_{j}
\circ^{1\mathbf{Pth}_{\boldsymbol{\mathcal{A}}^{(2)}}}_{s_{j}}
\mathfrak{P}^{(2)}_{j}
\right)
\right)
\right)_{j\in\bb{\mathbf{s}}}
\right)\right)
\tag{2}
\\&=
f_{s}\left(
\tau^{\mathbf{Pth}_{\boldsymbol{\mathcal{A}}^{(2)}}}
\left(\left(
\mathrm{ip}^{(2,X)@}_{s_{j}}\left(
\mathrm{CH}^{(2)}_{s_{j}}\left(
\mathfrak{Q}^{(2)}_{j}
\circ^{1\mathbf{Pth}_{\boldsymbol{\mathcal{A}}^{(2)}}}_{s_{j}}
\mathfrak{P}^{(2)}_{j}
\right)
\right)
\right)_{j\in\bb{\mathbf{s}}}
\right)\right)
\tag{3}
\\&=
\tau^{\mathbf{B}}
\left(\left(
f_{s_{j}}\left(
\mathrm{ip}^{(2,X)@}_{s_{j}}\left(
\mathrm{CH}^{(2)}_{s_{j}}\left(
\mathfrak{Q}^{(2)}_{j}
\circ^{1\mathbf{Pth}_{\boldsymbol{\mathcal{A}}^{(2)}}}_{s_{j}}
\mathfrak{P}^{(2)}_{j}
\right)
\right)
\right)
\right)_{j\in\bb{\mathbf{s}}}\right)
\tag{4}
\\&=
\tau^{\mathbf{B}}
\left(\left(
f_{s_{j}}\left(
\mathrm{ip}^{(2,X)@}_{s_{j}}\left(
\mathrm{CH}^{(2)}_{s_{j}}\left(
\mathfrak{Q}^{(2)}_{j}
\right)\right)\right)
\circ^{1\mathbf{B}}_{s_{j}}
\right.\right.
\\&\qquad\qquad\qquad\qquad\qquad\qquad\qquad\qquad
\left.\left.
f_{s_{j}}\left(
\mathrm{ip}^{(2,X)@}_{s_{j}}\left(
\mathrm{CH}^{(2)}_{s_{j}}\left(
\mathfrak{P}^{(2)}_{j}
\right)\right)\right)
\right)_{j\in\bb{\mathbf{s}}}\right)
\tag{5}
\\&=
\tau^{\mathbf{B}}
\left(\left(
f_{s_{j}}\left(
\mathrm{ip}^{(2,X)@}_{s_{j}}\left(
\mathrm{CH}^{(2)}_{s_{j}}\left(
\mathfrak{Q}^{(2)}_{j}
\right)\right)\right)
\right)_{j\in\bb{\mathbf{s}}}\right)
\circ^{1\mathbf{B}}_{s}
\\&\qquad\qquad\qquad\qquad\qquad\qquad\qquad
\tau^{\mathbf{B}}
\left(\left(
f_{s_{j}}\left(
\mathrm{ip}^{(2,X)@}_{s_{j}}\left(
\mathrm{CH}^{(2)}_{s_{j}}\left(
\mathfrak{P}^{(2)}_{j}
\right)\right)\right)
\right)_{j\in\bb{\mathbf{s}}}\right)
\tag{6}
\\&=
f_{s}\left(
\tau^{\mathbf{Pth}_{\boldsymbol{\mathcal{A}}^{(2)}}}
\left(\left(
\mathrm{ip}^{(2,X)@}_{s_{j}}\left(
\mathrm{CH}^{(2)}_{s_{j}}\left(
\mathfrak{Q}^{(2)}_{j}
\right)\right)
\right)_{j\in\bb{\mathbf{s}}}
\right)
\right)
\circ^{1\mathbf{B}}_{s}
\\&\qquad\qquad\qquad\qquad\qquad\qquad
f_{s}\left(
\tau^{\mathbf{Pth}_{\boldsymbol{\mathcal{A}}^{(2)}}}
\left(\left(
\mathrm{ip}^{(2,X)@}_{s_{j}}\left(
\mathrm{CH}^{(2)}_{s_{j}}\left(
\mathfrak{P}^{(2)}_{j}
\right)\right)
\right)_{j\in\bb{\mathbf{s}}}
\right)
\right)
\tag{7}
\\&=
f_{s}\left(
\tau^{
\mathbf{F}_{\Sigma^{\boldsymbol{\mathcal{A}}^{(2)}}}
(\mathbf{Pth}_{\boldsymbol{\mathcal{A}}^{(2)}})}
\left(\left(
\mathrm{ip}^{(2,X)@}_{s_{j}}\left(
\mathrm{CH}^{(2)}_{s_{j}}\left(
\mathfrak{Q}^{(2)}_{j}
\right)\right)
\right)_{j\in\bb{\mathbf{s}}}
\right)
\right)
\circ^{1\mathbf{B}}_{s}
\\&\qquad\qquad\qquad\qquad
f_{s}\left(
\tau^{
\mathbf{F}_{\Sigma^{\boldsymbol{\mathcal{A}}^{(2)}}}
(\mathbf{Pth}_{\boldsymbol{\mathcal{A}}^{(2)}})}
\left(\left(
\mathrm{ip}^{(2,X)@}_{s_{j}}\left(
\mathrm{CH}^{(2)}_{s_{j}}\left(
\mathfrak{P}^{(2)}_{j}
\right)\right)
\right)_{j\in\bb{\mathbf{s}}}
\right)
\right)
\tag{8}
\\&=
f_{s}\left(
\mathrm{ip}^{(2,X)@}_{s}\left(
\tau^{
\mathbf{T}_{\Sigma^{\boldsymbol{\mathcal{A}}^{(2)}}}(X)}
\left(\left(
\mathrm{CH}^{(2)}_{s_{j}}\left(
\mathfrak{Q}^{(2)}_{j}
\right)
\right)_{j\in\bb{\mathbf{s}}}
\right)
\right)
\right)
\circ^{1\mathbf{B}}_{s}
\\&\qquad\qquad\qquad\qquad\qquad
f_{s}\left(
\mathrm{ip}^{(2,X)@}_{s}\left(
\tau^{
\mathbf{T}_{\Sigma^{\boldsymbol{\mathcal{A}}^{(2)}}}(X)}
\left(\left(
\mathrm{CH}^{(2)}_{s_{j}}\left(
\mathfrak{P}^{(2)}_{j}
\right)
\right)_{j\in\bb{\mathbf{s}}}
\right)
\right)
\right)
\tag{9}
\\&=
f_{s}\left(
\mathrm{ip}^{(2,X)@}_{s}\left(
\mathrm{CH}^{(2)}_{s}\left(
\mathfrak{Q}^{(2)}\right)\right)\right)
\circ^{1\mathbf{B}}_{s}
f_{s}\left(
\mathrm{ip}^{(2,X)@}_{s}\left(
\mathrm{CH}^{(2)}_{s}\left(
\mathfrak{P}^{(2)}
\right)\right)
\right).
\tag{10}
\end{align*}
\end{flushleft}

In the just stated chain of equalities, the first equality follows from the previous discussion on the value of the second-order Curry-Howard mapping at $\mathfrak{Q}^{(2)}\circ^{1\mathbf{Pth}_{\boldsymbol{\mathcal{A}}^{(2)}}}_{s}\mathfrak{P}^{(2)}$; the second equality follow from the fact that $\mathrm{ip}^{(2,X)@}$ is a $\Sigma^{\boldsymbol{\mathcal{A}}^{(2)}}$-homomorphism, according to Definition~\ref{DDIp}; the third equality follows from the fact that, by Proposition~\ref{PDIpDCH}, for every $j\in\bb{\mathbf{s}}$ $\mathrm{ip}^{(2,X)@}_{s_{j}}(
\mathrm{CH}^{(2)}_{s_{j}}(
\mathfrak{P}^{(2)}_{j}
\circ^{1\mathbf{Pth}_{\boldsymbol{\mathcal{A}}^{(2)}}}_{s_{j}}
\mathfrak{P}^{(2)}_{j}
))$ is a second-order path. Thus, the interpretation of the operation symbol $\tau$ in the $\Sigma^{\boldsymbol{\mathcal{A}}^{(2)}}$-algebra $\mathbf{F}_{\Sigma^{\boldsymbol{\mathcal{A}}^{(2)}}}(\mathbf{Pth}_{\boldsymbol{\mathcal{A}}^{(2)}})$ becomes the corresponding interpretation of the $\tau$ operation symbol in the many-sorted partial $\Sigma^{\boldsymbol{\mathcal{A}}^{(2)}}$-algebra $\mathbf{Pth}_{\boldsymbol{\mathcal{A}}^{(2)}}$; the fourth equality follows from the fact that, by assumption $f$ is a $\Sigma^{\boldsymbol{\mathcal{A}}^{(2)}}$-homomorphism; the fifth equality follows by induction. Let us note that, for every $j\in\bb{\mathbf{s}}$, 
the pair $(\mathfrak{Q}^{(2)}_{j}\circ^{1\mathbf{Pth}_{\boldsymbol{\mathcal{A}}^{(2)}}}_{s_{j}}
\mathfrak{P}^{(2)}_{j},s_{j})$ $\prec_{\mathbf{Pth}_{\boldsymbol{\mathcal{A}}^{(2)}}}$-precedes $(\mathfrak{Q}^{(2)}\circ^{1\mathbf{Pth}_{\boldsymbol{\mathcal{A}}^{(2)}}}_{s}\mathfrak{P}^{(2)},s)$. Hence, we have that 
\begin{multline*}
f_{s_{j}}\left(
\mathrm{ip}^{(2,X)@}_{s_{j}}\left(
\mathrm{CH}^{(2)}_{s_{j}}\left(
\mathfrak{Q}^{(2)}_{j}
\circ^{1\mathbf{Pth}_{\boldsymbol{\mathcal{A}}^{(2)}}}_{s_{j}}
\mathfrak{P}^{(2)}_{j}
\right)\right)
\right)
\\=
f_{s_{j}}\left(
\mathrm{ip}^{(2,X)@}_{s_{j}}\left(
\mathrm{CH}^{(2)}_{s_{j}}\left(
\mathfrak{Q}^{(2)}_{j}
\right)\right)\right)
\circ^{1\mathbf{B}}_{s_{j}}
f_{s_{j}}\left(
\mathrm{ip}^{(2,X)@}_{s_{j}}\left(
\mathrm{CH}^{(2)}_{s_{j}}\left(
\mathfrak{P}^{(2)}_{j}
\right)\right)
\right);
\end{multline*}
the sixth equality follows from the fact that, since $\mathbf{B}$ is a partial many-sorted $\Sigma^{\boldsymbol{\mathcal{A}}^{(2)}}$-algebra in $\mathbf{PAlg}(\boldsymbol{\mathcal{E}}^{\boldsymbol{\mathcal{A}}^{(2)}})$, by either Axiom~\ref{DDVarB8} or Axiom~\ref{DDVarAB3}, the $1$-composition is compatible with the realization of the operation symbol $\tau$ in $\mathbf{B}$ ; the seventh equality follows from the fact that, by assumption, $f$ is a $\Sigma^{\boldsymbol{\mathcal{A}}^{(2)}}$-homomorphism; the eighth equality follows from the fact that,  by Proposition~\ref{PDIpDCH}, for every $j\in\bb{\mathbf{s}}$ $\mathrm{ip}^{(2,X)@}_{s_{j}}(
\mathrm{CH}^{(2)}_{s_{j}}(
\mathfrak{P}^{(2)}_{j}
\circ^{1\mathbf{Pth}_{\boldsymbol{\mathcal{A}}^{(2)}}}_{s_{j}}
\mathfrak{P}^{(2)}_{j}
))$ is a second-order path. Thus, the interpretation of the operation symbol $\tau$ in the $\Sigma^{\boldsymbol{\mathcal{A}}^{(2)}}$-algebra $\mathbf{F}_{\Sigma^{\boldsymbol{\mathcal{A}}^{(2)}}}(\mathbf{Pth}_{\boldsymbol{\mathcal{A}}^{(2)}})$ becomes the corresponding interpretation of the operation symbol $\tau$ in the many-sorted partial $\Sigma^{\boldsymbol{\mathcal{A}}^{(2)}}$-algebra $\mathbf{Pth}_{\boldsymbol{\mathcal{A}}^{(2)}}$; the ninth equality follows from the fact that $\mathrm{ip}^{(2,X)@}$ is a $\Sigma^{\boldsymbol{\mathcal{A}}^{(2)}}$-homomorphism, according to Definition~\ref{DDIp}; finally, the last equality recovers the value of the second-order Curry-Howard mapping at $\mathfrak{P}^{(2)}$ and $\mathfrak{Q}^{(2)}$, as we have discussed above.

The case~(2.3) follows.

It follows that Equation~\ref{PDVarKerE2} holds, i.e., the mapping $\mathrm{pr}^{
\equiv^{\llbracket 2\rrbracket}
}\circ\mathrm{ip}^{(2,X)@}\circ\mathrm{CH}^{(2)}$ is compatible with the $1$-composition.

All in all, we conclude that $\mathrm{pr}^{
\equiv^{\llbracket 2\rrbracket}
}\circ\mathrm{ip}^{(2,X)@}\circ\mathrm{CH}^{(2)}$ is a $\Sigma^{\boldsymbol{\mathcal{A}}^{(2)}}$-homomorphism.

It remains to prove that 
\[
\mathrm{Ker}(\mathrm{CH}^{(2)})\vee \Upsilon^{[1]}
\subseteq \mathrm{Ker}\left( \mathrm{pr}^{
\equiv^{\llbracket 2\rrbracket}
}\circ\mathrm{ip}^{(2,X)@}\circ\mathrm{CH}^{(2)}\right).
\]

In order to prove that $\mathrm{Ker}(\mathrm{CH}^{(2)})\vee \Upsilon^{[1]}$ is included in $\mathrm{Ker}( \mathrm{pr}^{
\equiv^{\llbracket 2\rrbracket}
}\circ\mathrm{ip}^{(2,X)@}\circ\mathrm{CH}^{(2)})$, it suffices to prove that both $\mathrm{Ker}(\mathrm{CH}^{(2)})$ and $\Upsilon^{[1]}$ are included in  $\mathrm{Ker}( \mathrm{pr}^{
\equiv^{\llbracket 2\rrbracket}
}\circ\mathrm{ip}^{(2,X)@}\circ\mathrm{CH}^{(2)})$.

\textsf{(i) $\mathrm{Ker}(\mathrm{CH}^{(2)})\subseteq \mathrm{Ker}( \mathrm{pr}^{
\equiv^{\llbracket 2\rrbracket}
}\circ\mathrm{ip}^{(2,X)@}\circ\mathrm{CH}^{(2)})$.}

Let $s$ be a sort in $S$ and let $\mathfrak{Q}^{(2)}$, $\mathfrak{P}^{(2)}$ be two second-order paths in $\mathrm{Pth}_{\boldsymbol{\mathcal{A}}^{(2)},s}$ satisfying that $(\mathfrak{Q}^{(2)},\mathfrak{P}^{(2)})\in\mathrm{Ker}(\mathrm{CH}^{(2)})_{s}$.

The following chain of equalities holds
\[
\mathrm{pr}^{
\equiv^{\llbracket 2\rrbracket}
}_{s}\left(\mathrm{ip}^{(2,X)@}_{s}\left(\mathrm{CH}^{(2)}_{s}\left(
\mathfrak{Q}^{(2)}
\right)\right)\right)
\\=
\mathrm{pr}^{
\equiv^{\llbracket 2\rrbracket}
}_{s}\left(\mathrm{ip}^{(2,X)@}_{s}\left(\mathrm{CH}^{(2)}_{s}\left(
\mathfrak{P}^{(2)}
\right)\right)\right).
\]

It follows directly from the fact that $\mathrm{CH}^{(2)}_{s}(\mathfrak{Q}^{(2)})=\mathrm{CH}^{(2)}_{s}(\mathfrak{P}^{(2)})$. Consequently, \[\mathrm{Ker}\left(\mathrm{CH}^{(2)}\right)\subseteq \mathrm{Ker}\left( \mathrm{pr}^{
\equiv^{\llbracket 2\rrbracket}
}\circ\mathrm{ip}^{(2,X)@}\circ\mathrm{CH}^{(2)}\right).\]

\textsf{(ii) $\Upsilon^{[1]}\subseteq \mathrm{Ker}( \mathrm{pr}^{
\equiv^{\llbracket 2\rrbracket}
}\circ\mathrm{ip}^{(2,X)@}\circ\mathrm{CH}^{(2)})$.}

Let us recall from Definition~\ref{DDUpsCong} that $\Upsilon^{[1]}$ is defined as the smallest $\Sigma^{\boldsymbol{\mathcal{A}}^{(2)}}$-congruence containing the relation $\Upsilon^{(1)}$, introduced in Definition~\ref{DDUps}. Therefore, since we have already proven that $\mathrm{pr}^{
\equiv^{\llbracket 2\rrbracket}
}\circ\mathrm{ip}^{(2,X)@}\circ\mathrm{CH}^{(2)}$ is a $\Sigma^{\boldsymbol{\mathcal{A}}^{(2)}}$-homomorphism, in order to check that $\Upsilon^{[1]}$ is included in $\mathrm{Ker}( \mathrm{pr}^{
\equiv^{\llbracket 2\rrbracket}
}\circ\mathrm{ip}^{(2,X)@}\circ\mathrm{CH}^{(2)})$ it suffices to prove that $\Upsilon^{(1)}$ is included in $\mathrm{Ker}( \mathrm{pr}^{
\equiv^{\llbracket 2\rrbracket}
}\circ\mathrm{ip}^{(2,X)@}\circ\mathrm{CH}^{(2)})$.

Following Definition~\ref{DDUps}, we will consider the different cases defining $\Upsilon^{(1)}$ individually.

\begin{enumerate}
\item For every sort $s\in S$ and every second-order path $\mathfrak{P}^{(2)}$ in $\mathrm{Pth}_{\boldsymbol{\mathcal{A}}^{(2)},s}$, 
\[
\left(
\mathfrak{P}^{(2)},
\mathfrak{P}^{(2)}\circ^{0\mathbf{Pth}_{\boldsymbol{\mathcal{A}}^{(2)}}}_{s}
\mathrm{sc}^{0\mathbf{Pth}_{\boldsymbol{\mathcal{A}}^{(2)}}}_{s}\left(
\mathfrak{P}^{(2)}
\right)
\right)
\in\Upsilon^{(1)}_{s}.
\]
\end{enumerate}

We want to check that 
\allowdisplaybreaks
\begin{multline*}
\mathrm{pr}^{\equiv^{\llbracket 2\rrbracket}}_{s}\left(
\mathrm{ip}^{(2,X)@}_{s}\left(
\mathrm{CH}^{(2)}_{s}\left(
\mathfrak{P}^{(2)}
\right)\right)\right)
\\=
\mathrm{pr}^{\equiv^{\llbracket 2\rrbracket}}_{s}\left(
\mathrm{ip}^{(2,X)@}_{s}\left(
\mathrm{CH}^{(2)}_{s}\left(
\mathfrak{P}^{(2)}\circ^{0\mathbf{Pth}_{\boldsymbol{\mathcal{A}}^{(2)}}}_{s}
\mathrm{sc}^{0\mathbf{Pth}_{\boldsymbol{\mathcal{A}}^{(2)}}}_{s}\left(
\mathfrak{P}^{(2)}
\right)
\right)\right)\right)
\end{multline*}

Developing both sides of the above equation, and taking into account Proposition~\ref{PDIpDCH}, we have that to check that the following equality between classes holds
\allowdisplaybreaks
\begin{multline*}
\left[
\mathrm{ip}^{(2,X)@}_{s}\left(
\mathrm{CH}^{(2)}_{s}\left(
\mathfrak{P}^{(2)}
\right)\right)
\right]_{\equiv^{\llbracket 2\rrbracket}}
\\=
\left[
\mathrm{ip}^{(2,X)@}_{s}\left(
\mathrm{CH}^{(2)}_{s}\left(
\mathfrak{P}^{(2)}\circ^{0\mathbf{Pth}_{\boldsymbol{\mathcal{A}}^{(2)}}}_{s}
\mathrm{sc}^{0\mathbf{Pth}_{\boldsymbol{\mathcal{A}}^{(2)}}}_{s}\left(
\mathfrak{P}^{(2)}
\right)
\right)\right)
\right]_{\equiv^{\llbracket 2\rrbracket}}.
\end{multline*}

But this is equivalent to prove that, for every many-sorted partial $\Sigma^{\boldsymbol{\mathcal{A}}^{(2)}}$-algebra $\mathbf{B}$ in $\mathbf{PAlg}(\boldsymbol{\mathcal{E}}^{\boldsymbol{\mathcal{A}}^{(2)}})$ and every $\Sigma^{\boldsymbol{\mathcal{A}}^{(2)}}$-homomorphism $f\colon\mathbf{Pth}_{\boldsymbol{\mathcal{A}}^{(2)}}\mor \mathbf{B}$ we have that 
\allowdisplaybreaks
\begin{multline*}
f^{\mathrm{Sch}}_{s}\left(
\mathrm{ip}^{(2,X)@}_{s}\left(
\mathrm{CH}^{(2)}_{s}\left(
\mathfrak{P}^{(2)}
\right)\right)
\right)
\\=
f^{\mathrm{Sch}}_{s}\left(
\mathrm{ip}^{(2,X)@}_{s}\left(
\mathrm{CH}^{(2)}_{s}\left(
\mathfrak{P}^{(2)}\circ^{0\mathbf{Pth}_{\boldsymbol{\mathcal{A}}^{(2)}}}_{s}
\mathrm{sc}^{0\mathbf{Pth}_{\boldsymbol{\mathcal{A}}^{(2)}}}_{s}\left(
\mathfrak{P}^{(2)}
\right)
\right)\right)
\right).
\end{multline*}

Let us recall, by Proposition~\ref{PDIpDCH}, that we are dealing with  proper second-order paths. Hence, the last mentioned equality reduces to
\allowdisplaybreaks
\begin{multline*}
f_{s}\left(
\mathrm{ip}^{(2,X)@}_{s}\left(
\mathrm{CH}^{(2)}_{s}\left(
\mathfrak{P}^{(2)}
\right)\right)
\right)
\\=
f_{s}\left(
\mathrm{ip}^{(2,X)@}_{s}\left(
\mathrm{CH}^{(2)}_{s}\left(
\mathfrak{P}^{(2)}\circ^{0\mathbf{Pth}_{\boldsymbol{\mathcal{A}}^{(2)}}}_{s}
\mathrm{sc}^{0\mathbf{Pth}_{\boldsymbol{\mathcal{A}}^{(2)}}}_{s}\left(
\mathfrak{P}^{(2)}
\right)
\right)\right)
\right).
\tag{E21}\label{PDVarKerE21}
\end{multline*}

Note that the following chain of equalities holds
\begin{flushleft}
$f_{s}\left(
\mathrm{ip}^{(2,X)@}_{s}\left(
\mathrm{CH}^{(2)}_{s}\left(
\mathfrak{P}^{(2)}\circ^{0\mathbf{Pth}_{\boldsymbol{\mathcal{A}}^{(2)}}}_{s}
\mathrm{sc}^{0\mathbf{Pth}_{\boldsymbol{\mathcal{A}}^{(2)}}}_{s}\left(
\mathfrak{P}^{(2)}
\right)
\right)\right)
\right)$
\allowdisplaybreaks
\begin{align*}
&=
f_{s}\left(
\mathrm{ip}^{(2,X)@}_{s}\left(
\mathrm{CH}^{(2)}_{s}\left(
\mathfrak{P}^{(2)}
\right)\right)
\circ^{0\mathbf{T}_{\boldsymbol{\mathcal{A}}^{(2)}}(X)}_{s}
\right.
\\&\qquad\qquad\qquad\qquad\qquad\qquad
\left.
\mathrm{sc}^{0\mathbf{T}_{\boldsymbol{\mathcal{A}}^{(2)}}(X)}_{s}\left(
\mathrm{ip}^{(2,X)@}_{s}\left(
\mathrm{CH}^{(2)}_{s}\left(
\mathfrak{P}^{(2)}
\right)\right)
\right)
\right)
\tag{1}
\\&=
f_{s}\left(
\mathrm{ip}^{(2,X)@}_{s}\left(
\mathrm{CH}^{(2)}_{s}\left(
\mathfrak{P}^{(2)}
\right)\right)
\right)
\circ^{0\mathbf{B}}_{s}
\mathrm{sc}^{0\mathbf{B}}_{s}\left(
f_{s}\left(
\mathrm{ip}^{(2,X)@}_{s}\left(
\mathrm{CH}^{(2)}_{s}\left(
\mathfrak{P}^{(2)}
\right)\right)
\right)
\right)
\tag{2}
\\&=
f_{s}\left(
\mathrm{ip}^{(2,X)@}_{s}\left(
\mathrm{CH}^{(2)}_{s}\left(
\mathfrak{P}^{(2)}
\right)\right)
\right).
\tag{3}
\end{align*}
\end{flushleft}

In the just stated chain of equalities, the first equality follows from Lemmas~\ref{LDIpDCHScZTgZ} and~\ref{LDIpDCHCompZ}; the second equality follows from the fact that, by assumption, $f$ is a $\Sigma^{\boldsymbol{\mathcal{A}}^{(2)}}$-homomorphism; finally, the last equality follows from the fact that $\mathbf{B}$ is a partial many-sorted $\Sigma^{\boldsymbol{\mathcal{A}}^{(2)}}$-algebra in $\mathbf{PAlg}(\boldsymbol{\mathcal{E}}^{\boldsymbol{\mathcal{A}}^{(2)}})$, then by  Axiom~\ref{DDVarA5}, the $0$-source is a neutral composition for the $0$-composition on the right.

This proves Case~(1).

\begin{enumerate}
\item[(2)] For every sort $s\in S$ and every second-order path $\mathfrak{P}^{(2)}$ in $\mathrm{Pth}_{\boldsymbol{\mathcal{A}}^{(2)},s}$, 
\[
\left(
\mathfrak{P}^{(2)},
\mathrm{tg}^{0\mathbf{Pth}_{\boldsymbol{\mathcal{A}}^{(2)}}}_{s}\left(
\mathfrak{P}^{(2)}
\right)
\circ^{0\mathbf{Pth}_{\boldsymbol{\mathcal{A}}^{(2)}}}_{s}
\mathfrak{P}^{(2)}
\right)
\in\Upsilon^{(1)}_{s}.
\]
\end{enumerate}

We want to check that 
\allowdisplaybreaks
\begin{multline*}
\mathrm{pr}^{\equiv^{\llbracket 2\rrbracket}}_{s}\left(
\mathrm{ip}^{(2,X)@}_{s}\left(
\mathrm{CH}^{(2)}_{s}\left(
\mathfrak{P}^{(2)}
\right)\right)\right)
\\=
\mathrm{pr}^{\equiv^{\llbracket 2\rrbracket}}_{s}\left(
\mathrm{ip}^{(2,X)@}_{s}\left(
\mathrm{CH}^{(2)}_{s}\left(
\mathrm{tg}^{0\mathbf{Pth}_{\boldsymbol{\mathcal{A}}^{(2)}}}_{s}\left(
\mathfrak{P}^{(2)}
\right)
\circ^{0\mathbf{Pth}_{\boldsymbol{\mathcal{A}}^{(2)}}}_{s}
\mathfrak{P}^{(2)}
\right)\right)\right)
\end{multline*}

This case follows by a similar argument to that presented for the Case~(1).

This proves Case~(2).

\begin{enumerate}
\item[(3)] For every sort $s\in S$ and every second-order paths $\mathfrak{P}^{(2)}$,$\mathfrak{Q}^{(2)}$ and $\mathfrak{R}^{(2)}$ in $\mathrm{Pth}_{\boldsymbol{\mathcal{A}}^{(2)},s}$, satisfying that
\begin{align*}
\mathrm{sc}^{(0,2)}_{s}\left(\mathfrak{R}^{(2)}\right)
&=
\mathrm{tg}^{(0,2)}_{s}\left(\mathfrak{Q}^{(2)}\right);
&
\mathrm{sc}^{(0,2)}_{s}\left(\mathfrak{Q}^{(2)}\right)
&=
\mathrm{tg}^{(0,2)}_{s}\left(\mathfrak{P}^{(2)}\right);
\end{align*}
then
\allowdisplaybreaks
\begin{multline*}
\left(
\mathfrak{R}^{(2)}
\circ^{0\mathbf{Pth}_{\boldsymbol{\mathcal{A}}^{(2)}}}_{s}
\left(
\mathfrak{Q}^{(2)}
\circ^{0\mathbf{Pth}_{\boldsymbol{\mathcal{A}}^{(2)}}}_{s}
\mathfrak{P}^{(2)}
\right),
\right.
\\
\left.
\left(
\mathfrak{R}^{(2)}
\circ^{0\mathbf{Pth}_{\boldsymbol{\mathcal{A}}^{(2)}}}_{s}
\mathfrak{Q}^{(2)}
\right)
\circ^{0\mathbf{Pth}_{\boldsymbol{\mathcal{A}}^{(2)}}}_{s}
\mathfrak{P}^{(2)}
\right)
\in\Upsilon^{(1)}_{s}.
\end{multline*}
\end{enumerate}

Note that the following chain of equalities holds
\begin{flushleft}
$f_{s}\left(
\mathrm{ip}^{(2,X)@}_{s}\left(
\mathrm{CH}^{(2)}_{s}\left(
\mathfrak{R}^{(2)}
\circ^{0\mathbf{Pth}_{\boldsymbol{\mathcal{A}}^{(2)}}}_{s}
\left(
\mathfrak{Q}^{(2)}
\circ^{0\mathbf{Pth}_{\boldsymbol{\mathcal{A}}^{(2)}}}_{s}
\mathfrak{P}^{(2)}
\right)
\right)\right)
\right)$
\allowdisplaybreaks
\begin{align*}
&=
f_{s}\left(
\mathrm{ip}^{(2,X)@}_{s}\left(
\mathrm{CH}^{(2)}_{s}\left(
\mathfrak{R}^{(2)}
\right)\right)
\circ^{0\mathbf{T}_{\Sigma^{\boldsymbol{\mathcal{A}}^{(2)}}}(X)}_{s}
\right.
\\&\qquad\qquad
\left.
\left(
\mathrm{ip}^{(2,X)@}_{s}\left(
\mathrm{CH}^{(2)}_{s}\left(
\mathfrak{Q}^{(2)}
\right)\right)
\circ^{0\mathbf{T}_{\Sigma^{\boldsymbol{\mathcal{A}}^{(2)}}}(X)}_{s}
\mathrm{ip}^{(2,X)@}_{s}\left(
\mathrm{CH}^{(2)}_{s}\left(
\mathfrak{P}^{(2)}
\right)\right)
\right)
\right)
\tag{1}
\\&=
f_{s}\left(
\mathrm{ip}^{(2,X)@}_{s}\left(
\mathrm{CH}^{(2)}_{s}\left(
\mathfrak{R}^{(2)}
\right)\right)\right)
\circ^{0\mathbf{B}}_{s}
\\&\qquad\qquad
\left(
f_{s}\left(
\mathrm{ip}^{(2,X)@}_{s}\left(
\mathrm{CH}^{(2)}_{s}\left(
\mathfrak{Q}^{(2)}
\right)\right)\right)
\circ^{0\mathbf{B}}_{s}
f_{s}\left(
\mathrm{ip}^{(2,X)@}_{s}\left(
\mathrm{CH}^{(2)}_{s}\left(
\mathfrak{P}^{(2)}
\right)\right)\right)
\right)
\tag{2}
\\&=
\left(
f_{s}\left(
\mathrm{ip}^{(2,X)@}_{s}\left(
\mathrm{CH}^{(2)}_{s}\left(
\mathfrak{R}^{(2)}
\right)\right)\right)
\circ^{0\mathbf{B}}_{s}
f_{s}\left(
\mathrm{ip}^{(2,X)@}_{s}\left(
\mathrm{CH}^{(2)}_{s}\left(
\mathfrak{Q}^{(2)}
\right)\right)\right)
\right)
\\&\qquad\qquad\qquad\qquad\qquad\qquad\qquad\qquad\qquad
\circ^{0\mathbf{B}}_{s}
f_{s}\left(
\mathrm{ip}^{(2,X)@}_{s}\left(
\mathrm{CH}^{(2)}_{s}\left(
\mathfrak{P}^{(2)}
\right)\right)\right)
\tag{3}
\\&=
f_{s}\left(
\left(
\mathrm{ip}^{(2,X)@}_{s}\left(
\mathrm{CH}^{(2)}_{s}\left(
\mathfrak{R}^{(2)}
\right)\right)
\circ^{0\mathbf{T}_{\Sigma^{\boldsymbol{\mathcal{A}}^{(2)}}}(X)}_{s}
\mathrm{ip}^{(2,X)@}_{s}\left(
\mathrm{CH}^{(2)}_{s}\left(
\mathfrak{Q}^{(2)}
\right)\right)
\right)
\right.
\\&\qquad\qquad\qquad\qquad\qquad\qquad\qquad\qquad
\left.
\circ^{0\mathbf{T}_{\Sigma^{\boldsymbol{\mathcal{A}}^{(2)}}}(X)}_{s}
\mathrm{ip}^{(2,X)@}_{s}\left(
\mathrm{CH}^{(2)}_{s}\left(
\mathfrak{P}^{(2)}
\right)\right)
\right)
\tag{4}
\\&=
f_{s}\left(
\mathrm{ip}^{(2,X)@}_{s}\left(
\mathrm{CH}^{(2)}_{s}\left(
\left(
\mathfrak{R}^{(2)}
\circ^{0\mathbf{Pth}_{\boldsymbol{\mathcal{A}}^{(2)}}}_{s}
\mathfrak{Q}^{(2)}
\right)
\circ^{0\mathbf{Pth}_{\boldsymbol{\mathcal{A}}^{(2)}}}_{s}
\mathfrak{P}^{(2)}
\right)
\right)
\right).
\tag{5}
\end{align*}
\end{flushleft}

In the just stated chain of equalities, the first equality follows from Lemma~\ref{LDIpDCHCompZ}; the second equality follows from the fact that, by assumption, $f$ is a $\Sigma^{\boldsymbol{\mathcal{A}}^{(2)}}$-homomorphism; the third equality follows from the fact that $\mathbf{B}$ is a partial many-sorted $\Sigma^{\boldsymbol{\mathcal{A}}^{(2)}}$-algebra in $\mathbf{PAlg}(\boldsymbol{\mathcal{E}}^{\boldsymbol{\mathcal{A}}^{(2)}})$, then by  Axiom~\ref{DDVarA6}, the $0$-composition is associative; the fourth equality follows from the fact that, by assumption, $f$ is a $\Sigma^{\boldsymbol{\mathcal{A}}^{(2)}}$-homomorphism;  finally, the last equality follows from Lemma~\ref{LDIpDCHCompZ}.

This proves Case~(3).

\begin{enumerate}
\item[(4)] For every word $\mathbf{s}$ and $s\in S$, for every operation symbol $\sigma\in \Sigma_{\mathbf{s},s}$, for every pair of families of second-order paths $(\mathfrak{P}^{(2)}_{j})_{j\in\bb{\mathbf{s}}}$ and $(\mathfrak{Q}^{(2)}_{j})_{j\in\bb{\mathbf{s}}}$ in $\mathrm{Pth}_{\boldsymbol{\mathcal{A}}^{(2)},\mathbf{s}}$ satisfying that, for every $j\in\bb{\mathbf{s}},$
\begin{align*}
\mathrm{sc}^{(0,2)}_{s_{j}}\left(\mathfrak{Q}^{(2)}_{j}\right)
&=
\mathrm{tg}^{(0,2)}_{s_{j}}\left(\mathfrak{P}^{(2)}_{j}\right);
\end{align*}
then
\allowdisplaybreaks
\begin{multline*}
\left(
\sigma^{\mathbf{Pth}_{\boldsymbol{\mathcal{A}}^{(2)}}}
\left(
\left(
\mathfrak{Q}^{(2)}_{j}
\circ^{0\mathbf{Pth}_{\boldsymbol{\mathcal{A}}^{(2)}}}_{s_{j}}
\mathfrak{P}^{(2)}_{j}
\right)_{j\in\bb{\mathbf{s}}}
\right)
\right.
\\
\left.
\sigma^{\mathbf{Pth}_{\boldsymbol{\mathcal{A}}^{(2)}}}
\left(
\left(
\mathfrak{Q}^{(2)}_{j}
\right)_{j\in\bb{\mathbf{s}}}
\right)
\circ^{0\mathbf{Pth}_{\boldsymbol{\mathcal{A}}^{(2)}}}_{s}
\sigma^{\mathbf{Pth}_{\boldsymbol{\mathcal{A}}^{(2)}}}
\left(
\left(
\mathfrak{P}^{(2)}_{j}
\right)_{j\in\bb{\mathbf{s}}}
\right)
\right)
\in\Upsilon^{(1)}_{s}.
\end{multline*}
\end{enumerate}

Note that the following chain of equalities holds
\begin{flushleft}
$f_{s}\left(
\mathrm{ip}^{(2,X)@}_{s}\left(
\mathrm{CH}^{(2)}_{s}\left(
\sigma^{\mathbf{Pth}_{\boldsymbol{\mathcal{A}}^{(2)}}}
\left(
\left(
\mathfrak{Q}^{(2)}_{j}
\circ^{0\mathbf{Pth}_{\boldsymbol{\mathcal{A}}^{(2)}}}_{s_{j}}
\mathfrak{P}^{(2)}_{j}
\right)_{j\in\bb{\mathbf{s}}}
\right)
\right)\right)
\right)$
\allowdisplaybreaks
\begin{align*}
&=
f_{s}\left(
\sigma^{\mathbf{T}_{\Sigma^{\boldsymbol{\mathcal{A}}^{(2)}}}(X)}\left(
\left(
\mathrm{ip}^{(2,X)@}_{s_{j}}\left(
\mathrm{CH}^{(2)}_{s_{j}}\left(
\mathfrak{Q}^{(2)}_{j}
\right)\right)
\right.\right.\right.
\\&\qquad\qquad\qquad\qquad\qquad\qquad
\left.\left.\left.
\circ^{0\mathbf{T}_{\Sigma^{\boldsymbol{\mathcal{A}}^{(2)}}}(X)}_{s_{j}}
\mathrm{ip}^{(2,X)@}_{s_{j}}\left(
\mathrm{CH}^{(2)}_{s_{j}}\left(
\mathfrak{P}^{(2)}_{j}
\right)\right)
\right)_{j\in\bb{\mathbf{s}}}
\right)
\right)
\tag{1}
\\&=
\sigma^{\mathbf{B}}\left(\left(
f_{s_{j}}\left(
\mathrm{ip}^{(2,X)@}_{s_{j}}\left(
\mathrm{CH}^{(2)}_{s_{j}}\left(
\mathfrak{Q}^{(2)}_{j}
\right)\right)
\right)
\right.\right.
\\&\qquad\qquad\qquad\qquad\qquad\qquad\qquad
\left.\left.
\circ^{0\mathbf{B}}_{s_{j}}
f_{s_{j}}\left(
\mathrm{ip}^{(2,X)@}_{s_{j}}\left(
\mathrm{CH}^{(2)}_{s_{j}}\left(
\mathfrak{P}^{(2)}_{j}
\right)\right)
\right)
\right)_{j\in\bb{\mathbf{s}}}
\right)
\tag{2}
\\&=
\sigma^{\mathbf{B}}\left(\left(
f_{s_{j}}\left(
\mathrm{ip}^{(2,X)@}_{s_{j}}\left(
\mathrm{CH}^{(2)}_{s_{j}}\left(
\mathfrak{Q}^{(2)}_{j}
\right)\right)
\right)
\right)_{j\in\bb{\mathbf{s}}}
\right)
\\&\qquad\qquad\qquad\qquad\qquad\quad
\circ^{0\mathbf{B}}_{s}
\sigma^{\mathbf{B}}\left(\left(
f_{s_{j}}\left(
\mathrm{ip}^{(2,X)@}_{s_{j}}\left(
\mathrm{CH}^{(2)}_{s_{j}}\left(
\mathfrak{P}^{(2)}_{j}
\right)\right)
\right)
\right)_{j\in\bb{\mathbf{s}}}
\right)
\tag{3}
\\&=
f_{s}\left(
\sigma^{\mathbf{T}_{\Sigma^{\boldsymbol{\mathcal{A}}^{(2)}}}(X)}
\left(\left(
\mathrm{ip}^{(2,X)@}_{s_{j}}\left(
\mathrm{CH}^{(2)}_{s_{j}}\left(
\mathfrak{Q}^{(2)}_{j}
\right)\right)
\right.\right.\right.
\\&\qquad\qquad\qquad\qquad\qquad\qquad
\left.\left.\left.
\circ^{0\mathbf{T}_{\Sigma^{\boldsymbol{\mathcal{A}}^{(2)}}}(X)}_{s_{j}}
\mathrm{ip}^{(2,X)@}_{s_{j}}\left(
\mathrm{CH}^{(2)}_{s_{j}}\left(
\mathfrak{P}^{(2)}_{j}
\right)\right)
\right)_{j\in\bb{\mathbf{s}}}\right)
\right)
\tag{4}
\\&=
f_{s}\left(
\mathrm{ip}^{(2,X)@}_{s}\left(
\mathrm{CH}^{(2)}_{s}\left(
\sigma^{\mathbf{Pth}_{\boldsymbol{\mathcal{A}}^{(2)}}}
\left(
\left(
\mathfrak{Q}^{(2)}_{j}
\right)_{j\in\bb{\mathbf{s}}}
\right)
\right.\right.\right.
\\&\qquad\qquad\qquad\qquad\qquad\qquad\qquad\qquad
\left.\left.\left.
\circ^{0\mathbf{Pth}_{\boldsymbol{\mathcal{A}}^{(2)}}}_{s}
\sigma^{\mathbf{Pth}_{\boldsymbol{\mathcal{A}}^{(2)}}}
\left(
\left(
\mathfrak{P}^{(2)}_{j}
\right)_{j\in\bb{\mathbf{s}}}
\right)
\right)
\right)
\right).
\tag{5}
\end{align*}
\end{flushleft}

In the just stated chain of equalities, the first equality follows from Lemmas~\ref{LDIpDCHSigma} and~\ref{LDIpDCHCompZ}; the second equality follows from the fact that, by assumption, $f$ is a $\Sigma^{\boldsymbol{\mathcal{A}}^{(2)}}$-homomorphism; the third equality follows from the fact that $\mathbf{B}$ is a partial many-sorted $\Sigma^{\boldsymbol{\mathcal{A}}^{(2)}}$-algebra in $\mathbf{PAlg}(\boldsymbol{\mathcal{E}}^{\boldsymbol{\mathcal{A}}^{(2)}})$, then by  Axiom~\ref{DDVarA8}, the $\sigma$ is a functor for the $0$-composition; the fourth equality follows from the fact that, by assumption, $f$ is a $\Sigma^{\boldsymbol{\mathcal{A}}^{(2)}}$-homomorphism;  finally, the last equality follows from Lemmas~\ref{LDIpDCHSigma} and~\ref{LDIpDCHCompZ}.

This proves Case~(4).

This completes the proof of Proposition~\ref{PDVarKer}.
\end{proof}

\begin{restatable}{remark}{RDQPUniv}
\label{RDQPUniv} Following Proposition~\ref{PDVarKer} and taking into account the Universal Property of the Quotient, there exists a unique $\Sigma^{\boldsymbol{\mathcal{A}}^{(2)}}$-homomorphism, that we will denote by $(\mathrm{pr}^{\equiv^{\llbracket 2\rrbracket}}
\circ
\mathrm{ip}^{(2,X)@}
\circ
\mathrm{CH}^{(2)})^{\natural}$, i.e.,
\[
(\mathrm{pr}^{\equiv^{\llbracket 2\rrbracket}}
\circ
\mathrm{ip}^{(2,X)@}
\circ
\mathrm{CH}^{(2)})^{\natural}
\colon 
\left\llbracket
\mathbf{Pth}_{\boldsymbol{\mathcal{A}}^{(2)}}
\right\rrbracket
\mor
\mathbf{T}_{\boldsymbol{\mathcal{E}}^{\boldsymbol{\mathcal{A}}^{(2)}}}\left(
\mathbf{Pth}_{\boldsymbol{\mathcal{A}}^{(2)}}
\right)
\]
satisfying that 
$
(\mathrm{pr}^{\equiv^{\llbracket 2\rrbracket}}
\circ
\mathrm{ip}^{(2,X)@}
\circ
\mathrm{CH}^{(2)})^{\natural}
\circ 
\mathrm{pr}^{\llbracket \cdot\rrbracket}
=
\mathrm{pr}^{\equiv^{\llbracket 2\rrbracket}}
\circ
\mathrm{ip}^{(2,X)@}
\circ
\mathrm{CH}^{(2)}.
$
\end{restatable}

\begin{restatable}{proposition}{PDQEta}
\label{PDQEta}
The  diagram of Figure~\ref{FDQEta} commutes, i.e., the following equality holds
$$
\mathrm{pr}^{\equiv^{\llbracket 2\rrbracket}}
\circ
\mathrm{ip}^{(2,X)@}
\circ
\mathrm{CH}^{(2)}
=
\eta^{(\llbracket 2 \rrbracket,\mathbf{Pth}_{\boldsymbol{\mathcal{A}}^{(2)}})}.
$$
\end{restatable}

\begin{figure}
\begin{tikzpicture}
[ACliment/.style={-{To [angle'=45, length=5.75pt, width=4pt, round]}
}, scale=0.8]
\node[] (P3) 		at 	(0,-4) 	[] 	{$\mathbf{Pth}_{\boldsymbol{\mathcal{A}}^{(2)}}$};
\node[] (TXQ) 	at 	(6,-4) 	[] 	{$\mathbf{T}_{\Sigma^{\boldsymbol{\mathcal{A}}^{(2)}}}
(X)$};
\node[] (SchP) 		at 	(6,-6) 	[] 	{$\mathbf{Sch}_{\boldsymbol{\mathcal{E}}^{\boldsymbol{\mathcal{A}}^{(2)}}}(\mathbf{Pth}_{\boldsymbol{\mathcal{A}}^{(2)}})$};
\node[] (TEP2) 	at 	(6,-8) 	[] 	{$\mathbf{T}_{\boldsymbol{\mathcal{E}}^{\boldsymbol{\mathcal{A}}^{(2)}}}
(\mathbf{Pth}_{\boldsymbol{\mathcal{A}}^{(2)}})$};

\draw[ACliment]  (TXQ) 	to node [right]	
{$\mathrm{ip}^{(2,X)@}$} (SchP);
\draw[ACliment]  (SchP) 	to node [right]	
{$\mathrm{pr}^{\equiv^{\llbracket 2\rrbracket}}$}  (TEP2);
\draw[ACliment, bend right=10]  (P3) 	to node [below left]	
{$\eta^{(\llbracket 2 \rrbracket,\mathbf{Pth}_{\boldsymbol{\mathcal{A}}^{(2)}})}$} (TEP2);

\draw[ACliment]  (P3) 	to node [above]	
{$\mathrm{CH}^{(2)}$} (TXQ);
\end{tikzpicture}
\caption{The maps in Proposition~\ref{PDQEta}.}
\label{FDQEta}
\end{figure}

\begin{proof}
Let $s$ be a sort in $S$ and $\mathfrak{P}^{(2)}$ a second-order path in $\mathrm{Pth}_{\boldsymbol{\mathcal{A}}^{(2)},s}$. We need to prove that
$$
\left[\mathfrak{P}^{(2)}\right]_{\equiv^{\llbracket 2\rrbracket}_{s}}=
\left[
\mathrm{ip}^{(2,X)@}_{s}\left(
\mathrm{CH}^{(2)}_{s}\left(
\mathfrak{P}^{(2)}
\right)\right)\right]_{\equiv^{\llbracket 2\rrbracket}_{s}}.
$$
But this is equivalent to prove that, for every many-sorted partial $\Sigma^{\boldsymbol{\mathcal{A}}^{(2)}}$-algebra $\mathbf{B}$ in $\mathbf{PAlg}(\boldsymbol{\mathcal{E}}^{\boldsymbol{\mathcal{A}}^{(2)}})$ and every $\Sigma^{\boldsymbol{\mathcal{A}}^{(2)}}$-homomorphism $f\colon\mathbf{Pth}_{\boldsymbol{\mathcal{A}}^{(2)}}\mor\mathbf{B}$ we have
$$
f^{\mathrm{Sch}}_{s}\left(
\mathfrak{P}^{(2)}
\right)=
f^{\mathrm{Sch}}_{s}\left(
\mathrm{ip}^{(2,X)@}_{s}\left(
\mathrm{CH}^{(2)}_{s}\left(
\mathfrak{P}^{(2)}
\right)\right)\right).
$$

However, by Proposition~\ref{PDIpDCH}, the element $\mathrm{ip}^{(2,X)@}_{s}(\mathrm{CH}^{(2)}_{s}(\mathfrak{P}^{(2)}))$ is a second-order path. Therefore, the last equation reduces to the equation
$$
f_{s}\left(
\mathfrak{P}^{(2)}\right)
=f_{s}\left(
\mathrm{ip}^{(2,X)@}_{s}\left(
\mathrm{CH}^{(2)}_{s}\left(
\mathfrak{P}^{(2)}
\right)\right)\right).
$$

We prove this statement by Artinian induction on $(\coprod\mathrm{Pth}_{\boldsymbol{\mathcal{A}}^{(2)}},\leq_{\mathbf{Pth}_{\boldsymbol{\mathcal{A}}^{(2)}}})$.

\textsf{Base step of the Artinian induction}.

Let $(\mathfrak{P}^{(2)},s)$ be a minimal element in $(\coprod\mathrm{Pth}_{\boldsymbol{\mathcal{A}}^{(2)}}, \leq_{\mathbf{Pth}_{\boldsymbol{\mathcal{A}}^{(2)}}})$. Then, by Proposition~\ref{PDMinimal}, the second-order path $\mathfrak{P}^{(2)}$ is either a $(2,[1])$-identity second-order path on a simple term, or a second-order echelon. In any case, by Corollaries\ref{CDIpDUId} or~\ref{CDIpDCHOneStep}, if $\mathfrak{P}^{(2)}$ is a $(2,[1])$-identity second-order path or a one-step second-order path, then the following equation holds
$$
\mathfrak{P}^{(2)}=\mathrm{ip}^{(2,X)@}_{s}\left(
\mathrm{CH}^{(2)}_{s}\left(
\mathfrak{P}^{(2)}\right)\right).
$$
Therefore, the statement holds for second-order paths of length $0$ or $1$.

This case follows easily.

\textsf{Inductive step of the Artinian induction}.

Let $(\mathfrak{P}^{(2)},s)$ be a non-minimal element in $(\coprod\mathrm{Pth}_{\boldsymbol{\mathcal{A}}^{(2)}}, \leq_{\mathbf{Pth}_{\boldsymbol{\mathcal{A}}^{(2)}}})$. Let us suppose that, for every sort $t\in S$ and every second-order path $\mathfrak{Q}^{(2)}\in\mathrm{Pth}_{\boldsymbol{\mathcal{A}}^{(2)},t}$, if $(\mathfrak{Q}^{(2)},t)$ $<_{\mathbf{Pth}_{\boldsymbol{\mathcal{A}}^{(2)}}}$-precedes $(\mathfrak{P}^{(2)},s)$, then the statement holds for $\mathfrak{Q}^{(2)}$, i.e., we have that
$$
f_{t}\left(
\mathfrak{Q}^{(2)}
\right)=f_{t}\left(
\mathrm{ip}^{(2,X)@}_{t}\left(
\mathrm{CH}^{(2)}_{t}\left(
\mathfrak{Q}^{(2)}
\right)\right)\right).
$$

Let us recall that, by the discussion presented in the base case, we can also assume that $\mathfrak{P}^{(2)}$ is not a $(2,1)$-identity second-order path. Since $(\mathfrak{P}^{(2)},s)$ is a non-minimal element in $(\coprod\mathrm{Pth}_{\boldsymbol{\mathcal{A}}^{(2)}}, \leq_{\mathbf{Pth}_{\boldsymbol{\mathcal{A}}}})$, we have, by Lemma~\ref{LDOrdI}, that $\mathfrak{P}^{(2)}$ is either~(1) a second-order path of length strictly greater than one containing at least one second-order echelon or~(2) an echelonless second-order path.

If~(1), then let $i\in \bb{\mathfrak{P}^{(2)}}$ be the first index for which the one-step subpath $\mathfrak{P}^{(2),i,i}$ of $\mathfrak{P}^{(2)}$ is a second-order echelon. We distinguish the cases (1.1) $i=0$ and (1.2) $i>0$.

If~(1.1), i.e., if $\mathfrak{P}^{(2)}$ is a second-order path of length strictly greater than one having a second-order echelon on its first step  then, according to Definition~\ref{DDCH}, the value of the second-order Curry-Howard mapping at $\mathfrak{P}^{(2)}$ is given by
$$
\mathrm{CH}^{(2)}_{s}\left(
\mathfrak{P}^{(2)}
\right)=
\mathrm{CH}^{(2)}_{s}\left(
\mathfrak{P}^{(2),1,\bb{\mathfrak{P}^{(2)}}-1}
\right)
\circ_{s}^{1\mathbf{T}_{\Sigma^{\boldsymbol{\mathcal{A}}^{(2)}}}(X)}
\mathrm{CH}^{(2)}_{s}\left(
\mathfrak{P}^{(2),0,0}
\right).
$$

Thus, the following chain of equalities holds
\begin{flushleft}
$f_{s}\left(
\mathrm{ip}^{(2,X)@}_{s}\left(
\mathrm{CH}^{(2)}_{s}\left(
\mathfrak{P}^{(2)}
\right)\right)\right)$
\allowdisplaybreaks
\begin{align*}
&=
f_{s}\left(
\mathrm{ip}^{(2,X)@}_{s}\left(
\mathrm{CH}^{(2)}_{s}\left(
\mathfrak{P}^{(2),1,\bb{\mathfrak{P}^{(2)}}-1}
\right)
\circ_{s}^{1\mathbf{T}_{\Sigma^{\boldsymbol{\mathcal{A}}^{(2)}}}(X)}
\mathrm{CH}^{(2)}_{s}\left(
\mathfrak{P}^{(2),0,0}
\right)
\right)\right)
\tag{1}
\\&=
f_{s}\left(
\mathrm{ip}^{(2,X)@}_{s}\left(
\mathrm{CH}^{(2)}_{s}\left(
\mathfrak{P}^{(2),1,\bb{\mathfrak{P}^{(2)}}-1}
\right)\right)
\circ_{s}^{1\mathbf{F}_{\Sigma^{\boldsymbol{\mathcal{A}}^{(2)}}}
(\mathbf{Pth}_{\boldsymbol{\mathcal{A}}^{(2)}})}
\right.
\\&\qquad\qquad\qquad\qquad\qquad\qquad\qquad\qquad\qquad\qquad
\left.
\mathrm{ip}^{(2,X)@}_{s}\left(
\mathrm{CH}^{(2)}_{s}\left(
\mathfrak{P}^{(2),0,0}
\right)\right)\right)
\tag{2}
\\&=
f_{s}\left(
\mathrm{ip}^{(2,X)@}_{s}\left(
\mathrm{CH}^{(2)}_{s}\left(
\mathfrak{P}^{(2),1,\bb{\mathfrak{P}^{(2)}}-1}
\right)\right)
\circ_{s}^{1\mathbf{Pth}_{\boldsymbol{\mathcal{A}}^{(2)}}}
\right.
\\&\qquad\qquad\qquad\qquad\qquad\qquad\qquad\qquad\qquad\qquad
\left.
\mathrm{ip}^{(2,X)@}_{s}\left(
\mathrm{CH}^{(2)}_{s}\left(
\mathfrak{P}^{(2),0,0}
\right)\right)\right)
\tag{3}
\\&=
f_{s}\left(
\mathrm{ip}^{(2,X)@}_{s}\left(
\mathrm{CH}^{(2)}_{s}\left(
\mathfrak{P}^{(2),1,\bb{\mathfrak{P}^{(2)}}-1}
\right)\right)\right)
\circ_{s}^{1\mathbf{B}}
\\&\qquad\qquad\qquad\qquad\qquad\qquad\qquad\qquad\qquad\quad
f_{s}\left(
\mathrm{ip}^{(2,X)@}_{s}\left(
\mathrm{CH}^{(2)}_{s}\left(
\mathfrak{P}^{(2),0,0}
\right)\right)\right)
\tag{4}
\\&=
f_{s}\left(
\mathfrak{P}^{(2),1,\bb{\mathfrak{P}^{(2)}}-1}
\right)
\circ_{s}^{1\mathbf{B}}
f_{s}\left(
\mathfrak{P}^{(2),0,0}
\right)
\tag{5}
\\&=
f_{s}\left(
\mathfrak{P}^{(2),1,\bb{\mathfrak{P}^{(2)}}-1}
\circ_{s}^{1\mathbf{\mathbf{Pth}_{\boldsymbol{\mathcal{A}}^{(2)}}}}
\mathfrak{P}^{(2),0,0}
\right)
\tag{6}
\\&=
f_{s}\left(
\mathfrak{P}^{(2)}
\right).
\tag{7}
\end{align*}
\end{flushleft}

In the just stated chain of equalities, the first equality holds because it consists, simply, in unravelling the definition of the second-order Curry-Howard mapping at $\mathfrak{P}^{(2)}$; the second equality holds because $\mathrm{ip}^{(2,X)@}$ is a $\Sigma^{\boldsymbol{\mathcal{A}}^{(2)}}$-homomorphism, according to Definition~\ref{DDIp}; the third equation follows from the fact that, in virtue of Lemma~\ref{PDIpDCH},  both $\mathrm{ip}^{(2,X)@}_{s}(
\mathrm{CH}^{(2)}_{s}(\mathfrak{P}^{(2),1,\bb{\mathfrak{P}^{(2)}}-1}))$ and $\mathrm{ip}^{(2,X)@}_{s}(
\mathrm{CH}^{(2)}_{s}(\mathfrak{P}^{(2),0,0}))$ are second-order paths. Thus, the interpretation of the operation symbol for the $1$-composition, $\circ^{1}_{s}$,  in the $\Sigma^{\boldsymbol{\mathcal{A}}^{(2)}}$-algebra $\mathbf{F}_{\Sigma^{\boldsymbol{\mathcal{A}}^{(2)}}}(\mathbf{Pth}_{\boldsymbol{\mathcal{A}}^{(2)}})$ becomes the corresponding interpretation of the $1$-composition $\circ^{1}_{s}$ in the $\Sigma^{\boldsymbol{\mathcal{A}}^{(2)}}$-algebra $\mathbf{Pth}_{\boldsymbol{\mathcal{A}}^{(2)}}$; the fourth equality holds because $f$ is a $\Sigma^{\boldsymbol{\mathcal{A}}^{(2)}}$-homomorphism;
the fifth equality holds by induction, since the pairs $(\mathfrak{P}^{(2),1,\bb{\mathfrak{P}^{(2)}}-1},s)$ and $(\mathfrak{P}^{(2),0,0},s)$ $\prec_{\mathbf{Pth}_{\boldsymbol{\mathcal{A}}}}$-precede $(\mathfrak{P}^{(2)},s)$; the fifth equation holds since $f$ is a $\Sigma^{\boldsymbol{\mathcal{A}}^{(2)}}$-homomorphism; and the last equality holds because in it we, simply, are recovering the original second-order path $\mathfrak{P}^{(2)}$ as the composition of its $1$-constituents.

If~(1.2), i.e., if $\mathfrak{P}^{(2)}$ is a second-order path of length strictly greater than one having its first second-order echelon on a step $i\in\bb{\mathfrak{P}}$ different from the initial one, then the value of the second-order Curry-Howard mapping at $\mathfrak{P}^{(2)}$ is given by
$$
\mathrm{CH}^{(2)}_{s}\left(
\mathfrak{P}^{(2)}
\right)=
\mathrm{CH}^{(2)}_{s}\left(
\mathfrak{P}^{(2),i,\bb{\mathfrak{P}^{(2)}}-1}
\right)
\circ_{s}^{1\mathbf{T}_{\Sigma^{\boldsymbol{\mathcal{A}}^{(2)}}}(X)}
\mathrm{CH}^{(2)}_{s}\left(
\mathfrak{P}^{(2),0,i-1}
\right).
$$

Thus, the following chain of equalities holds
\begin{flushleft}
$f_{s}\left(
\mathrm{ip}^{(2,X)@}_{s}\left(
\mathrm{CH}^{(2)}_{s}\left(
\mathfrak{P}^{(2)}
\right)\right)\right)$
\allowdisplaybreaks
\begin{align*}
&=
f_{s}\left(
\mathrm{ip}^{(2,X)@}_{s}\left(
\mathrm{CH}^{(2)}_{s}\left(
\mathfrak{P}^{(2),i,\bb{\mathfrak{P}^{(2)}}-1}
\right)
\circ_{s}^{1\mathbf{T}_{\Sigma^{\boldsymbol{\mathcal{A}}^{(2)}}}(X)}
\mathrm{CH}^{(2)}_{s}\left(
\mathfrak{P}^{(2),0,i-1}
\right)\right)\right)
\tag{1}
\\&=
f_{s}\left(
\mathrm{ip}^{(2,X)@}_{s}\left(
\mathrm{CH}^{(2)}_{s}\left(
\mathfrak{P}^{(2),i,\bb{\mathfrak{P}^{(2)}}-1}
\right)\right)
\circ_{s}^{1\mathbf{F}_{\Sigma^{\boldsymbol{\mathcal{A}}^{(2)}}}
(\mathbf{Pth}_{\boldsymbol{\mathcal{A}}^{(2)}})}
\right.
\\&\qquad\qquad\qquad\qquad\qquad\qquad\qquad\qquad\qquad\qquad
\left.
\mathrm{ip}^{(2,X)@}_{s}\left(
\mathrm{CH}^{(2)}_{s}\left(
\mathfrak{P}^{(2),0,i-1}
\right)\right)\right)
\tag{2}
\\&=
f_{s}\left(
\mathrm{ip}^{(2,X)@}_{s}\left(
\mathrm{CH}^{(2)}_{s}\left(
\mathfrak{P}^{(2),i,\bb{\mathfrak{P}^{(2)}}-1}
\right)\right)
\circ_{s}^{1\mathbf{Pth}_{\boldsymbol{\mathcal{A}}^{(2)}}}
\right.
\\&\qquad\qquad\qquad\qquad\qquad\qquad\qquad\qquad\qquad\qquad
\left.
\mathrm{ip}^{(2,X)@}_{s}\left(
\mathrm{CH}^{(2)}_{s}\left(
\mathfrak{P}^{(2),0,i-1}
\right)\right)\right)
\tag{3}
\\&=
f_{s}\left(
\mathrm{ip}^{(2,X)@}_{s}\left(
\mathrm{CH}^{(2)}_{s}\left(
\mathfrak{P}^{(2),i,\bb{\mathfrak{P}^{(2)}}-1}
\right)\right)\right)
\circ_{s}^{1\mathbf{B}}
\\&\qquad\qquad\qquad\qquad\qquad\qquad\qquad\qquad\qquad
f_{s}\left(
\mathrm{ip}^{(2,X)@}_{s}\left(
\mathrm{CH}^{(2)}_{s}\left(
\mathfrak{P}^{(2),0,i-1}
\right)\right)\right)
\tag{4}
\\&=
f_{s}\left(
\mathfrak{P}^{(2),i,\bb{\mathfrak{P}^{(2)}}-1}
\right)
\circ_{s}^{1\mathbf{B}}
f_{s}\left(
\mathfrak{P}^{(2),0,i-1}
\right)
\tag{5}
\\&=
f_{s}\left(
\mathfrak{P}^{(2),i,\bb{\mathfrak{P}^{(2)}}-1}
\circ_{s}^{1\mathbf{\mathbf{Pth}_{\boldsymbol{\mathcal{A}}^{(2)}}}}
\mathfrak{P}^{(2),0,i-1}
\right)
\tag{6}
\\&=
f_{s}\left(
\mathfrak{P}^{(2)}
\right).
\tag{7}
\end{align*}
\end{flushleft}

In the just stated chain of equalities, the first equality holds because it consists, simply, in unravelling the definition of the second-order Curry-Howard mapping at $\mathfrak{P}^{(2)}$; the second equality holds because $\mathrm{ip}^{(2,X)@}$ is a $\Sigma^{\boldsymbol{\mathcal{A}}^{(2)}}$-homomorphism, according to Definition~\ref{DDIp}; the third equation follows from the fact that, in virtue of Lemma~\ref{PDIpDCH},  both $\mathrm{ip}^{(2,X)@}_{s}(
\mathrm{CH}^{(2)}_{s}(\mathfrak{P}^{(2),i,\bb{\mathfrak{P}^{(2)}}-1}))$ and $\mathrm{ip}^{(2,X)@}_{s}(
\mathrm{CH}^{(2)}_{s}(\mathfrak{P}^{(2),0,i-1}))$ are proper second-order paths. Thus, the interpretation of the operation symbol for the $1$-composition, $\circ^{1}_{s}$,  in the $\Sigma^{\boldsymbol{\mathcal{A}}^{(2)}}$-algebra $\mathbf{F}_{\Sigma^{\boldsymbol{\mathcal{A}}^{(2)}}}(\mathbf{Pth}_{\boldsymbol{\mathcal{A}}^{(2)}})$ becomes the corresponding interpretation of the $1$-composition $\circ^{1}_{s}$ in the $\Sigma^{\boldsymbol{\mathcal{A}}^{(2)}}$-algebra $\mathbf{Pth}_{\boldsymbol{\mathcal{A}}^{(2)}}$; the fourth equality holds because $f$ is a $\Sigma^{\boldsymbol{\mathcal{A}}^{(2)}}$-homomorphism;
the fifth equality holds by induction, since the pairs $(\mathfrak{P}^{(2),i,\bb{\mathfrak{P}^{(2)}}-1},s)$ and $(\mathfrak{P}^{(2),0,i-1},s)$ $\prec_{\mathbf{Pth}_{\boldsymbol{\mathcal{A}}}}$-precede $(\mathfrak{P}^{(2)},s)$; the fifth equation holds since $f$ is a $\Sigma^{\boldsymbol{\mathcal{A}}^{(2)}}$-homomorphism; finally, the last equality holds because in it we, simply, are recovering the original second-order path $\mathfrak{P}^{(2)}$ as the composition of its $1$-constituents.

Case~(1) follows.

If~(2), i.e., if $\mathfrak{P}^{(2)}$ is an echelonless second-order path, it could be the case that (2.1) $\mathfrak{P}^{(2)}$ is an echelonless second-order path that is not head-constant, or (2.2) $\mathfrak{P}^{(2)}$ is a head-constant non-coherent echelonless second-order path or (2.3) $\mathfrak{P}^{(2)}$ is a head-constant coherent echelonless second-order path.

If~(2.1), let $i\in\bb{\mathfrak{P}^{(2)}}$ be the greatest index for which $\mathfrak{P}^{(2),0,i}$ is a head-constant second-order path. Then, according to Definition~\ref{DDCH}, we have that
$$
\mathrm{CH}^{(2)}_{s}\left(
\mathfrak{P}^{(2)}
\right)=
\mathrm{CH}^{(2)}_{s}\left(
\mathfrak{P}^{(2),i+1,\bb{\mathfrak{P}^{(2)}}-1}
\right)
\circ_{s}^{1\mathbf{T}_{\Sigma^{\boldsymbol{\mathcal{A}}^{(2)}}}(X)}
\mathrm{CH}^{(2)}_{s}\left(
\mathfrak{P}^{(2),0,i}
\right).
$$

Thus, the following chain of equalities holds
\begin{flushleft}
$f_{s}\left(
\mathrm{ip}^{(2,X)@}_{s}\left(
\mathrm{CH}^{(2)}_{s}\left(
\mathfrak{P}^{(2)}
\right)\right)\right)$
\allowdisplaybreaks
\begin{align*}
&=
f_{s}\left(
\mathrm{ip}^{(2,X)@}_{s}\left(
\mathrm{CH}^{(2)}_{s}\left(
\mathfrak{P}^{(2),i+1,\bb{\mathfrak{P}^{(2)}}-1}
\right)
\circ_{s}^{1\mathbf{T}_{\Sigma^{\boldsymbol{\mathcal{A}}^{(2)}}}(X)}
\mathrm{CH}^{(2)}_{s}\left(
\mathfrak{P}^{(2),0,i}
\right)\right)\right)
\tag{1}
\\&=
f_{s}\left(
\mathrm{ip}^{(2,X)@}_{s}\left(
\mathrm{CH}^{(2)}_{s}\left(
\mathfrak{P}^{(2),i+1,\bb{\mathfrak{P}^{(2)}}-1}
\right)\right)
\circ_{s}^{1\mathbf{F}_{\Sigma^{\boldsymbol{\mathcal{A}}^{(2)}}}
(\mathbf{Pth}_{\boldsymbol{\mathcal{A}}^{(2)}})}
\right.
\\&\qquad\qquad\qquad\qquad\qquad\qquad\qquad\qquad\qquad\qquad
\left.
\mathrm{ip}^{(2,X)@}_{s}\left(
\mathrm{CH}^{(2)}_{s}\left(
\mathfrak{P}^{(2),0,i}
\right)\right)\right)
\tag{2}
\\&=
f_{s}\left(
\mathrm{ip}^{(2,X)@}_{s}\left(
\mathrm{CH}^{(2)}_{s}\left(
\mathfrak{P}^{(2),i+1,\bb{\mathfrak{P}^{(2)}}-1}
\right)\right)
\circ_{s}^{1\mathbf{Pth}_{\boldsymbol{\mathcal{A}}^{(2)}}}
\right.
\\&\qquad\qquad\qquad\qquad\qquad\qquad\qquad\qquad\qquad\qquad
\left.
\mathrm{ip}^{(2,X)@}_{s}\left(
\mathrm{CH}^{(2)}_{s}\left(
\mathfrak{P}^{(2),0,i}
\right)\right)\right)
\tag{3}
\\&=
f_{s}\left(
\mathrm{ip}^{(2,X)@}_{s}\left(
\mathrm{CH}^{(2)}_{s}\left(
\mathfrak{P}^{(2),i+1,\bb{\mathfrak{P}^{(2)}}-1}
\right)\right)\right)
\circ_{s}^{1\mathbf{B}}
\\&\qquad\qquad\qquad\qquad\qquad\qquad\qquad\qquad\qquad
f_{s}\left(
\mathrm{ip}^{(2,X)@}_{s}\left(
\mathrm{CH}^{(2)}_{s}\left(
\mathfrak{P}^{(2),0,i}
\right)\right)\right)
\tag{4}
\\&=
f_{s}\left(
\mathfrak{P}^{(2),i+1,\bb{\mathfrak{P}^{(2)}}-1}
\right)
\circ_{s}^{1\mathbf{B}}
f_{s}\left(
\mathfrak{P}^{(2),0,i}
\right)
\tag{5}
\\&=
f_{s}\left(
\mathfrak{P}^{(2),i+1,\bb{\mathfrak{P}^{(2)}}-1}
\circ_{s}^{1\mathbf{\mathbf{Pth}_{\boldsymbol{\mathcal{A}}^{(2)}}}}
\mathfrak{P}^{(2),0,i}
\right)
\tag{6}
\\&=
f_{s}\left(
\mathfrak{P}^{(2)}
\right).
\tag{7}
\end{align*}
\end{flushleft}

In the just stated chain of equalities, the first equality holds because it consists, simply, in unravelling the definition of the second-order Curry-Howard mapping at $\mathfrak{P}^{(2)}$; the second equality holds because $\mathrm{ip}^{(2,X)@}$ is a $\Sigma^{\boldsymbol{\mathcal{A}}^{(2)}}$-homomorphism, according to Definition~\ref{DDIp}; the third equation follows from the fact that, in virtue of Lemma~\ref{PDIpDCH},  both $\mathrm{ip}^{(2,X)@}_{s}(
\mathrm{CH}^{(2)}_{s}(\mathfrak{P}^{(2),i+1,\bb{\mathfrak{P}^{(2)}}-1}))$ and $\mathrm{ip}^{(2,X)@}_{s}(
\mathrm{CH}^{(2)}_{s}(\mathfrak{P}^{(2),0,i}))$ are proper second-order paths. Thus, the interpretation of the operation symbol for the $1$-composition, $\circ^{1}_{s}$,  in the $\Sigma^{\boldsymbol{\mathcal{A}}^{(2)}}$-algebra $\mathbf{F}_{\Sigma^{\boldsymbol{\mathcal{A}}^{(2)}}}(\mathbf{Pth}_{\boldsymbol{\mathcal{A}}^{(2)}})$ becomes the corresponding interpretation of the $1$-composition $\circ^{1}_{s}$ in the $\Sigma^{\boldsymbol{\mathcal{A}}^{(2)}}$-algebra $\mathbf{Pth}_{\boldsymbol{\mathcal{A}}^{(2)}}$; the fourth equality holds because $f$ is a $\Sigma^{\boldsymbol{\mathcal{A}}^{(2)}}$-homomorphism;
the fifth equality holds by induction, since the pairs $(\mathfrak{P}^{(2),i+1,\bb{\mathfrak{P}^{(2)}}-1},s)$ and $(\mathfrak{P}^{(2),0,i},s)$ $\prec_{\mathbf{Pth}_{\boldsymbol{\mathcal{A}}}}$-precede $(\mathfrak{P}^{(2)},s)$; the fifth equation holds since $f$ is a $\Sigma^{\boldsymbol{\mathcal{A}}^{(2)}}$-homomorphism; finally, the last equality holds because in it we, simply, are recovering the original second-order path $\mathfrak{P}^{(2)}$ as the composition of its $1$-constituents.

The case of $\mathfrak{P}^{(2)}$ being an echelonless second-order path that is not head-constant follows.

If~(2.2), let $i\in\bb{\mathfrak{P}^{(2)}}$ be the greatest index for which $\mathfrak{P}^{(2),0,i}$ is a coherent head-constant second-order path. Then, according to Definition~\ref{DDCH}, we have that
$$
\mathrm{CH}^{(2)}_{s}\left(
\mathfrak{P}^{(2)}
\right)=
\mathrm{CH}^{(2)}_{s}\left(
\mathfrak{P}^{(2),i+1,\bb{\mathfrak{P}^{(2)}}-1}
\right)
\circ_{s}^{1\mathbf{T}_{\Sigma^{\boldsymbol{\mathcal{A}}^{(2)}}}(X)}
\mathrm{CH}^{(2)}_{s}\left(
\mathfrak{P}^{(2),0,i}
\right).
$$

Thus, the following chain of equalities holds
\begin{flushleft}
$f_{s}\left(
\mathrm{ip}^{(2,X)@}_{s}\left(
\mathrm{CH}^{(2)}_{s}\left(
\mathfrak{P}^{(2)}
\right)\right)\right)$
\allowdisplaybreaks
\begin{align*}
&=
f_{s}\left(
\mathrm{ip}^{(2,X)@}_{s}\left(
\mathrm{CH}^{(2)}_{s}\left(
\mathfrak{P}^{(2),i+1,\bb{\mathfrak{P}^{(2)}}-1}
\right)
\circ_{s}^{1\mathbf{T}_{\Sigma^{\boldsymbol{\mathcal{A}}^{(2)}}}(X)}
\mathrm{CH}^{(2)}_{s}\left(
\mathfrak{P}^{(2),0,i}
\right)\right)\right)
\tag{1}
\\&=
f_{s}\left(
\mathrm{ip}^{(2,X)@}_{s}\left(
\mathrm{CH}^{(2)}_{s}\left(
\mathfrak{P}^{(2),i+1,\bb{\mathfrak{P}^{(2)}}-1}
\right)\right)
\circ_{s}^{1\mathbf{F}_{\Sigma^{\boldsymbol{\mathcal{A}}^{(2)}}}
(\mathbf{Pth}_{\boldsymbol{\mathcal{A}}^{(2)}})}
\right.
\\&\qquad\qquad\qquad\qquad\qquad\qquad\qquad\qquad\qquad\qquad
\left.
\mathrm{ip}^{(2,X)@}_{s}\left(
\mathrm{CH}^{(2)}_{s}\left(
\mathfrak{P}^{(2),0,i}
\right)\right)\right)
\tag{2}
\\&=
f_{s}\left(
\mathrm{ip}^{(2,X)@}_{s}\left(
\mathrm{CH}^{(2)}_{s}\left(
\mathfrak{P}^{(2),i+1,\bb{\mathfrak{P}^{(2)}}-1}
\right)\right)
\circ_{s}^{1\mathbf{Pth}_{\boldsymbol{\mathcal{A}}^{(2)}}}
\right.
\\&\qquad\qquad\qquad\qquad\qquad\qquad\qquad\qquad\qquad\qquad
\left.
\mathrm{ip}^{(2,X)@}_{s}\left(
\mathrm{CH}^{(2)}_{s}\left(
\mathfrak{P}^{(2),0,i}
\right)\right)\right)
\tag{3}
\\&=
f_{s}\left(
\mathrm{ip}^{(2,X)@}_{s}\left(
\mathrm{CH}^{(2)}_{s}\left(
\mathfrak{P}^{(2),i+1,\bb{\mathfrak{P}^{(2)}}-1}
\right)\right)\right)
\circ_{s}^{1\mathbf{B}}
\\&\qquad\qquad\qquad\qquad\qquad\qquad\qquad\qquad\qquad
f_{s}\left(
\mathrm{ip}^{(2,X)@}_{s}\left(
\mathrm{CH}^{(2)}_{s}\left(
\mathfrak{P}^{(2),0,i}
\right)\right)\right)
\tag{4}
\\&=
f_{s}\left(
\mathfrak{P}^{(2),i+1,\bb{\mathfrak{P}^{(2)}}-1}
\right)
\circ_{s}^{1\mathbf{B}}
f_{s}\left(
\mathfrak{P}^{(2),0,i}
\right)
\tag{5}
\\&=
f_{s}\left(
\mathfrak{P}^{(2),i+1,\bb{\mathfrak{P}^{(2)}}-1}
\circ_{s}^{1\mathbf{\mathbf{Pth}_{\boldsymbol{\mathcal{A}}^{(2)}}}}
\mathfrak{P}^{(2),0,i}
\right)
\tag{6}
\\&=
f_{s}\left(
\mathfrak{P}^{(2)}
\right).
\tag{7}
\end{align*}
\end{flushleft}

In the just stated chain of equalities, the first equality holds because it consists, simply, in unravelling the definition of the second-order Curry-Howard mapping at $\mathfrak{P}^{(2)}$; the second equality holds because $\mathrm{ip}^{(2,X)@}$ is a $\Sigma^{\boldsymbol{\mathcal{A}}^{(2)}}$-homomorphism, according to Definition~\ref{DDIp}; the third equation follows from the fact that, in virtue of Lemma~\ref{PDIpDCH},  both $\mathrm{ip}^{(2,X)@}_{s}(
\mathrm{CH}^{(2)}_{s}(\mathfrak{P}^{(2),i+1,\bb{\mathfrak{P}^{(2)}}-1}))$ and $\mathrm{ip}^{(2,X)@}_{s}(
\mathrm{CH}^{(2)}_{s}(\mathfrak{P}^{(2),0,i}))$ are proper second-order paths. Thus, the interpretation of the operation symbol for the $1$-composition, $\circ^{1}_{s}$,  in the $\Sigma^{\boldsymbol{\mathcal{A}}^{(2)}}$-algebra $\mathbf{F}_{\Sigma^{\boldsymbol{\mathcal{A}}^{(2)}}}(\mathbf{Pth}_{\boldsymbol{\mathcal{A}}^{(2)}})$ becomes the corresponding interpretation of the $1$-composition $\circ^{1}_{s}$ in the $\Sigma^{\boldsymbol{\mathcal{A}}^{(2)}}$-algebra $\mathbf{Pth}_{\boldsymbol{\mathcal{A}}^{(2)}}$; the fourth equality holds because $f$ is a $\Sigma^{\boldsymbol{\mathcal{A}}^{(2)}}$-homomorphism;
the fifth equality holds by induction, since the pairs $(\mathfrak{P}^{(2),i+1,\bb{\mathfrak{P}^{(2)}}-1},s)$ and $(\mathfrak{P}^{(2),0,i},s)$ $\prec_{\mathbf{Pth}_{\boldsymbol{\mathcal{A}}}}$-precede $(\mathfrak{P}^{(2)},s)$; the fifth equation holds since $f$ is a $\Sigma^{\boldsymbol{\mathcal{A}}^{(2)}}$-homomorphism; finally, the last equality holds because in it we, simply, are recovering the original second-order path $\mathfrak{P}^{(2)}$ as the composition of its $1$-constituents.

The case of $\mathfrak{P}^{(2)}$ being a head-constant echelonless second-order path that is not coherent follows.

If~(2.3), then there exists a unique word $\mathbf{s}\in S^{\star}-\{\lambda\}$ and a unique operation symbol $\tau\in\Sigma^{\boldsymbol{\mathcal{A}}}_{\mathbf{s},s}$ associated to $\mathfrak{P}^{(2)}$. Let $(\mathfrak{P}^{(2)}_{j})_{j\in\bb{\mathbf{s}}}$ be the family of second-order paths in $\mathbf{Pth}_{\boldsymbol{\mathcal{A}},\mathbf{s}}$ which, in virtue of Lemma~\ref{LDPthExtract}, we can extract from $\mathfrak{P}^{(2)}$. Then, according to Definition~\ref{DDCH}, the image of the second-order Curry-Howard mapping at $\mathfrak{P}^{(2)}$ is given by
$$
\mathrm{CH}^{(2)}_{s}\left(
\mathfrak{P}^{(2)}
\right)
=
\tau^{\mathbf{T}_{\Sigma^{\boldsymbol{\mathcal{A}}^{(2)}}}(X)}
\left(\left(
\mathrm{CH}^{(2)}_{s_{j}}\left(
\mathfrak{P}^{(2)}_{j}
\right)\right)_{j\in\bb{\mathbf{s}}}
\right).
$$

Thus, the following chain of equalities holds
\begin{flushleft}
$f_{s}\left(
\mathrm{ip}^{(2,X)@}_{s}\left(
\mathrm{CH}^{(2)}_{s}\left(
\mathfrak{P}^{(2)}
\right)\right)\right)$
\allowdisplaybreaks
\begin{align*}
\qquad
&=
f_{s}\left(
\mathrm{ip}^{(2,X)@}_{s}\left(
\tau^{\mathbf{T}_{\Sigma^{\boldsymbol{\mathcal{A}}^{(2)}}}(X)}
\left(\left(
\mathrm{CH}^{(2)}_{s_{j}}\left(
\mathfrak{P}^{(2)}_{j}
\right)\right)_{j\in\bb{\mathbf{s}}}
\right)\right)\right)
\tag{1}
\\&=
f_{s}\left(
\tau^{\mathbf{F}_{\Sigma^{\boldsymbol{\mathcal{A}}^{(2)}}}(
\mathbf{Pth}_{\boldsymbol{\mathcal{A}}^{(2)}})}
\left(\left(
\mathrm{ip}^{(2,X)@}_{s_{j}}\left(
\mathrm{CH}^{(2)}_{s_{j}}\left(
\mathfrak{P}^{(2)}_{j}
\right)\right)\right)_{j\in\bb{\mathbf{s}}}
\right)\right)
\tag{2}
\\&=
f_{s}\left(
\tau^{\mathbf{Pth}_{\boldsymbol{\mathcal{A}}^{(2)}}}
\left(\left(
\mathrm{ip}^{(2,X)@}_{s_{j}}\left(
\mathrm{CH}^{(2)}_{s_{j}}\left(
\mathfrak{P}^{(2)}_{j}
\right)\right)\right)_{j\in\bb{\mathbf{s}}}
\right)\right)
\tag{3}
\\&=
\tau^{\mathbf{B}}
\left(\left(
f_{s_{j}}\left(
\mathrm{ip}^{(2,X)@}_{s_{j}}\left(
\mathrm{CH}^{(2)}_{s_{j}}\left(
\mathfrak{P}^{(2)}_{j}
\right)\right)\right)\right)_{j\in\bb{\mathbf{s}}}
\right)
\tag{4}
\\&=
\tau^{\mathbf{B}}
\left(\left(
f_{s_{j}}\left(
\mathfrak{P}^{(2)}_{j}
\right)\right)_{j\in\bb{\mathbf{s}}}
\right)
\tag{5}
\\&=
f_{s}\left(
\tau^{\mathbf{Pth}_{\boldsymbol{\mathcal{A}}^{(2)}}}
\left(\left(
\mathfrak{P}^{(2)}_{j}
\right)_{j\in\bb{\mathbf{s}}}
\right)\right).
\tag{6}
\end{align*}
\end{flushleft}

In the just stated chain of equalities, the first equality holds because it consists, simply, in unravelling the definition of the second-order Curry-Howard mapping at $\mathfrak{P}^{(2)}$; the second equality holds because $\mathrm{ip}^{(2,X)@}$ is a $\Sigma^{\boldsymbol{\mathcal{A}}^{(2)}}$-homomorphism, according to Definition~\ref{DDIp}; the third equation follows from the fact that, in virtue of Lemma~\ref{PDIpDCH},  for every $j\in\bb{\mathbf{s}}$, the element
$\mathrm{ip}^{(2,X)@}_{s_{j}}(
\mathrm{CH}^{(2)}_{s_{j}}(\mathfrak{P}^{(2)}_{j}
))$  is a second-order path. Thus, the interpretation of the operation symbol $\tau$,  in the $\Sigma^{\boldsymbol{\mathcal{A}}^{(2)}}$-algebra $\mathbf{F}_{\Sigma^{\boldsymbol{\mathcal{A}}^{(2)}}}(\mathbf{Pth}_{\boldsymbol{\mathcal{A}}^{(2)}})$ becomes the corresponding interpretation of $\tau$ in the $\Sigma^{\boldsymbol{\mathcal{A}}^{(2)}}$-algebra $\mathbf{Pth}_{\boldsymbol{\mathcal{A}}^{(2)}}$; the fourth equality holds because $f$ is a $\Sigma^{\boldsymbol{\mathcal{A}}^{(2)}}$-homomorphism; the fifth equality holds by induction, since, for every $j\in\bb{\mathbf{s}}$, the pair $((\mathfrak{P}^{(2)}_{j},s_{j}),(\mathfrak{P}^{(2)},s))$ is in $\prec_{\mathbf{Pth}_{\boldsymbol{\mathcal{A}}}}$; finally, the last equality holds because $f$  is a $\Sigma^{\boldsymbol{\mathcal{A}}^{(2)}}$-homomorphism.

It remains to prove that, for a coherent head-constant echelonless path $\mathfrak{P}^{(2)}$, the following equality holds
$$
f_{s}\left(
\mathfrak{P}^{(2)}
\right)=f_{s}\left(
\tau^{\mathbf{Pth}_{\boldsymbol{\mathcal{A}}^{(2)}}}
\left(\left(
\mathfrak{P}^{(2)}_{j}
\right)_{j\in\bb{\mathbf{s}}}
\right)\right).
$$

Let $\mathbf{c}$ be the unique word of $S^{\star}$ for which $\mathfrak{P}^{(2)}$ is a second-order path in $\mathrm{Pth}_{\mathbf{c},\boldsymbol{\mathcal{A}}^{(2)},s}$. 
By Proposition~\ref{PDPthRecons}, we have that $\mathfrak{P}^{(2)}$ can be univocally represented as a $1$-composition of one-step second-order paths as
\begin{align*}
\mathfrak{P}^{(2)}&=
\mathfrak{P}^{(2),\bb{\mathbf{c}}-1,\bb{\mathbf{c}}-1}
\circ^{1\mathbf{Pth}_{\boldsymbol{\mathcal{A}}^{(2)}}}_{s}
\cdots
\circ^{1\mathbf{Pth}_{\boldsymbol{\mathcal{A}}^{(2)}}}_{s}
\mathfrak{P}^{(2),0,0}.
\tag{E1}
\label{PDQEtaE1}
\end{align*}

Moreover, since $\mathfrak{P}^{(2)}$ is a coherent head-constant echelonless second-order path, we have, by Lemma~\ref{LDPthExtract}, that there exists a partition of $\mathbf{c}$ into a family $(\mathbf{c}_{j})_{j\in\bb{\mathbf{s}}}$ in $(S^{\star})^{\bb{\mathbf{s}}}$ such that, for every $j\in\bb{\mathbf{s}}$, $\mathfrak{P}^{(2)}_{j}$ is a second-order path in $\mathrm{Pth}_{\mathbf{c}_{j},\boldsymbol{\mathcal{A}}^{(2)},s_{j}}$.

Let us note that, for every $i\in \bb{\mathbf{c}}$, the one-step second-order subpath $\mathfrak{P}^{(2),i,i}$ of $\mathfrak{P}^{(2)}$ is a coherent head-constant echelonless second-order path associated to the operation symbol $\tau$. Now, recovering the terminology used in Lemma~\ref{LDPthExtract}, the index $i$ must have type $(j_{i},k_{i})$, for some $j_{i}\in\bb{\mathbf{s}}$ and some $k_{i}\in\bb{\mathbf{c}_{j_{i}}}$, meaning that the unique first-order translation appearing on the one-step second-order path $\mathfrak{P}^{(2),i,i}$ has its derived first-order translation occurring at position $j_{i}$, this fact being the $k_{i}$-th time that occurs, i.e.,  there exists a family of paths $(\mathfrak{P}_{i,l})_{l\in j_{i}}\in \prod_{l\in j_{i}}\mathrm{Pth}_{\boldsymbol{\mathcal{A}},s_{l}}$ and a family of paths $(\mathfrak{P}_{i,l})_{l\in \bb{\mathbf{s}}-(j_{i}+1)}\in \prod_{l\in \bb{\mathbf{s}}-(j_{i}+1)}\mathrm{Pth}_{\boldsymbol{\mathcal{A}},s_{l}}$ and a first-order translation $T'_{i}\in \mathrm{Tl}_{c_{i}}(\mathbf{T}_{\Sigma^{\boldsymbol{\mathcal{A}}}}(X))_{s_{j_{i}}}$ such that
\begin{multline*}
T_{i}=\tau^{\mathbf{T}_{\Sigma^{\boldsymbol{\mathcal{A}}}}(X)}\left(
\mathrm{CH}^{(1)}_{s_{0}}\left(\mathfrak{P}_{i,0}\right),
\cdots,
\mathrm{CH}^{(1)}_{s_{j_{i}-1}}\left(\mathfrak{P}_{i,j_{i}-1}\right),
\right.
\\
\left.
T'_{i},
\mathrm{CH}^{(1)}_{s_{j_{i}+1}}\left(\mathfrak{P}_{i,j_{i}+1}\right),
\cdots,
\mathrm{CH}^{(1)}_{s_{\bb{\mathbf{s}}-1}}\left(\mathfrak{P}_{i,\bb{\mathbf{s}}-1}\right)
\right).
\end{multline*}

Let $((\mathfrak{P}^{(2),i,i})_{j})_{j\in\bb{\mathbf{s}}}$ be the family of second-order paths that we can extract from $\mathfrak{P}^{(2),i,i}$ in virtue of Lemma~\ref{LDPthExtract}. What we obtain, following the proof of Lemma~\ref{LDPthExtract}, is the following family of second-order paths
\allowdisplaybreaks
\begin{align*}
\left(\mathfrak{P}^{(2),i,i}\right)_{0}&=
\mathrm{ip}^{(2,[1])\sharp}_{s_{0}}\left(
\left[
\mathrm{CH}^{(1)}_{s_{0}}\left(
\mathfrak{P}_{i,0}
\right)
\right]_{s_{0}}
\right);
\\
\vdots&\qquad\vdots
\\
\left(\mathfrak{P}^{(2),i,i}\right)_{j_{i}-1}&=
\mathrm{ip}^{(2,[1])\sharp}_{s_{j_{i}-1}}\left(
\left[
\mathrm{CH}^{(1)}_{s_{j_{i}-1}}\left(
\mathfrak{P}_{i,j_{i}-1}
\right)
\right]_{s_{j_{i}-1}}
\right);
\\
\left(\mathfrak{P}^{(2),i,i}\right)_{j_{i}}&=
\left(\mathfrak{P}^{(2)}_{j_{i}}\right)^{k_{i},k_{i}};
\\
\left(\mathfrak{P}^{(2),i,i}\right)_{j_{i}+1}&=
\mathrm{ip}^{(2,[1])\sharp}_{s_{j_{i}+1}}\left(
\left[
\mathrm{CH}^{(1)}_{s_{j_{i}+1}}\left(
\mathfrak{P}_{i,j_{i}+1}
\right)
\right]_{s_{j_{i}+1}}
\right);
\\
\vdots&\qquad\vdots
\\
\left(\mathfrak{P}^{(2),i,i}\right)_{\bb{\mathbf{s}}-1}&=
\mathrm{ip}^{(2,[1])\sharp}_{s_{\bb{\mathbf{s}}-1}}\left(
\left[
\mathrm{CH}^{(1)}_{s_{\bb{\mathbf{s}}-1}}\left(
\mathfrak{P}_{i,\bb{\mathbf{s}}-1}
\right)
\right]_{s_{\bb{\mathbf{s}}-1}}
\right).
\end{align*}

That is, we recover $(2,[1])$-identity second-order paths on every index in $\bb{\mathbf{s}}$ different from $j_{i}$, whilst on index $j_{i}$ we recover, precisely, the $k_{i}$-th step of the $j_{i}$-th second-order path that we can extract from $\mathfrak{P}^{(2)}$, i.e., $(\mathfrak{P}^{(2)}_{j_{i}})^{k_{i},k_{i}}$.

Two questions regarding the families of second-order paths $((\mathfrak{P}^{(2),i,i}_{j})_{j\in \bb{\mathbf{s}}})_{i\in \bb{\mathbf{c}}}$ arise. On the first hand, let us note that, since $\mathfrak{P}^{(2),i,i}$ is a coherent head-constant one-step echelonless second-order path, it follows, by Corollary~\ref{CDUStep}, that
\begin{align*}
\mathfrak{P}^{(2),i,i}
&=
\tau^{\mathbf{Pth}_{\boldsymbol{\mathcal{A}}^{(2)}}}
\left(
\left(
\left(\mathfrak{P}^{(2),i,i}\right)_{j}
\right)_{j\in\bb{\mathbf{s}}}
\right).
\tag{E2}
\label{PDQEtaE2}
\end{align*}

Moreover, following Proposition~\ref{PDPthRecons}, we can reconstruct, for every $j\in\bb{\mathbf{s}}$, the second-order path $\mathfrak{P}^{(2)}_{j}$, as the $1$-composition of its one-step subpaths as follows
$$
\mathfrak{P}^{(2)}_{j}
=
\left(\mathfrak{P}^{(2)}_{j}\right)^{\bb{\mathbf{c}_{j}}-1, \bb{\mathbf{c}_{j}}-1}
\circ^{1\mathbf{Pth}_{\boldsymbol{\mathcal{A}}^{(2)}}}_{s_{j}}
\cdots
\circ^{1\mathbf{Pth}_{\boldsymbol{\mathcal{A}}^{(2)}}}_{s_{j}}
\left(\mathfrak{P}^{(2)}_{j}\right)^{0,0}.
$$

The following equation is also a valid description of $\mathfrak{P}^{(2)}_{j}$ as a  $1$-composition
\begin{align*}
\mathfrak{P}^{(2)}_{j}&=
\left(\mathfrak{P}^{(2),\bb{\mathbf{c}}-1,\bb{\mathbf{c}}-1}\right)_{j}
\circ^{1\mathbf{Pth}_{\boldsymbol{\mathcal{A}}^{(2)}}}_{s_{j}}
\cdots
\circ^{1\mathbf{Pth}_{\boldsymbol{\mathcal{A}}^{(2)}}}_{s_{j}}
\left(\mathfrak{P}^{(2),0,0}\right)_{j}.
\tag{E3}
\label{PDQEtaE3}
\end{align*}

In this last equation we are $1$-composing all the second-order paths at position $j$ that we can extract from the one-step subpaths $\mathfrak{P}^{(2),i,i}$, for every $i\in\bb{\mathbf{c}}$. Note that, for those indexes $i\in\bb{\mathbf{c}}$, that are not of type $j$, the second-order path $(\mathfrak{P}^{(2),i,i})_{j}$ will simply be a $(2,[1])$-identity second-order path that acts as a neutral element for the $1$-composition.

All in all, the following chain of equalities holds
\begin{flushleft}
$f_{s}\left(
\mathfrak{P}^{(2)}
\right)$
\allowdisplaybreaks
\begin{align*}
&=
f_{s}\left(
\mathfrak{P}^{(2),\bb{\mathbf{c}}-1,\bb{\mathbf{c}}-1}
\circ^{1\mathbf{Pth}_{\boldsymbol{\mathcal{A}}^{(2)}}}_{s}
\cdots
\circ^{1\mathbf{Pth}_{\boldsymbol{\mathcal{A}}^{(2)}}}_{s}
\mathfrak{P}^{(2),0,0}
\right)
\tag{1}
\\&=
f_{s}\left(
\mathfrak{P}^{(2),\bb{\mathbf{c}}-1,\bb{\mathbf{c}}-1}
\right)
\circ^{1\mathbf{B}}_{s}
\cdots
\circ^{1\mathbf{B}}_{s}
f_{s}
\Big(
\mathfrak{P}^{(2),0,0}
\Big)
\tag{2}
\\&=
f_{s}\left(
\tau^{\mathbf{Pth}_{\boldsymbol{\mathcal{A}}^{(2)}}}
\left(
\left(
\left(
\mathfrak{P}^{(2),\bb{\mathbf{c}}-1,\bb{\mathbf{c}}-1}
\right)_{j}
\right)_{j\in\bb{\mathbf{s}}}
\right)\right)
\circ^{1\mathbf{B}}_{s}
\cdots
\circ^{1\mathbf{B}}_{s}
\\&
\qquad\qquad\qquad\qquad\qquad\qquad\qquad\qquad
f_{s}\left(
\tau^{\mathbf{Pth}_{\boldsymbol{\mathcal{A}}^{(2)}}}
\left(
\left(
\left(
\mathfrak{P}^{(2),0,0}
\right)_{j}
\right)_{j\in\bb{\mathbf{s}}}
\right)\right)
\tag{3}
\\&=
\tau^{\mathbf{B}}
\left(
\left(
f_{s_{j}}\left(
\left(
\mathfrak{P}^{(2),\bb{\mathbf{c}}-1,\bb{\mathbf{c}}-1}
\right)_{j}
\right)
\right)_{j\in\bb{\mathbf{s}}}\right)
\circ^{1\mathbf{B}}_{s}
\cdots
\circ^{1\mathbf{B}}_{s}
\\&
\qquad\qquad\qquad\qquad\qquad\qquad\qquad\qquad\qquad
\tau^{\mathbf{B}}
\left(
\left(
f_{s_{j}}\left(
\left(
\mathfrak{P}^{(2),0,0}
\right)_{j}
\right)
\right)_{j\in\bb{\mathbf{s}}}\right)
\tag{4}
\\&=
\tau^{\mathbf{B}}
\left(\left(
f_{s_{j}}\left(
\left(
\mathfrak{P}^{(2),\bb{\mathbf{c}}-1,\bb{\mathbf{c}}-1}
\right)_{j}
\right)
\circ^{1\mathbf{B}}_{s_{j}}
\cdots 
\circ^{1\mathbf{B}}_{s_{j}}
f_{s_{j}}\left(
\left(
\mathfrak{P}^{(2),0,0}
\right)_{j}
\right)
\right)_{j\in\bb{\mathbf{s}}}
\right)
\tag{5}
\\
&=
\tau^{\mathbf{B}}
\left(\left(
f_{s_{j}}\left(
\left(\mathfrak{P}^{(2),\bb{\mathbf{c}}-1,\bb{\mathbf{c}}-1}\right)_{j}
\circ^{1\mathbf{Pth}_{\boldsymbol{\mathcal{A}}^{(2)}}}_{s_{j}}
\cdots
\circ^{1\mathbf{Pth}_{\boldsymbol{\mathcal{A}}^{(2)}}}_{s_{j}}
\left(\mathfrak{P}^{(2),0,0}\right)_{j}
\right)
\right)_{j\in\bb{\mathbf{s}}}
\right)
\tag{6}
\\
&=
\tau^{\mathbf{B}}
\left(\left(
f_{s_{j}}\left(
\mathfrak{P}^{(2)}_{j}
\right)\right)_{j\in\bb{\mathbf{s}}}
\right)
\tag{7}
\\&=
f_{s}\left(\tau^{\mathbf{Pth}_{\boldsymbol{\mathcal{A}}^{(2)}}}
\left(\left(
\mathfrak{P}^{(2)}_{j}
\right)_{j\in\bb{\mathbf{s}}}
\right)\right).
\tag{8}
\end{align*}
\end{flushleft}

In the just stated chain of equalities, the first equality holds because $\mathfrak{P}^{(2)}$ can be decomposed into one-step second-order paths, as shown in Equation~\ref{PDQEtaE1}; the second equality holds because $f$ is a $\Sigma^{\boldsymbol{\mathcal{A}}^{(2)}}$-homomorphism; the third equation follows
from Equation~\ref{PDQEtaE2}; the fourth equality holds because $f$ is a $\Sigma^{\boldsymbol{\mathcal{A}}^{(2)}}$-homomorphism; the fifth equality holds because the many-sorted partial $\Sigma^{\boldsymbol{\mathcal{A}}^{(2)}}$-algebra $\mathbf{B}$ is in $\mathbf{PAlg}(\boldsymbol{\mathcal{E}}^{\boldsymbol{\mathcal{A}}^{(2)}})$ thus, by either Axiom~\ref{DDVarB8} or Axiom~\ref{DDVarAB3}, the $1$-composition,  $\circ^{1}_{s}$, is compatible with the interpretation of the operation $\tau$ in $\mathbf{B}$, this operation being an operation $\sigma\in \Sigma_{\mathbf{s},s}$ or the $0$-composition; the sixth equality holds because $f$ is a $\Sigma^{\boldsymbol{\mathcal{A}}^{(2)}}$-homomorphism; the seventh equality holds according to Equation~\ref{PDQEtaE3}; finally, the last equality holds because $f$ is a $\Sigma^{\boldsymbol{\mathcal{A}}^{(2)}}$-homomorphism.

This completes Case~(2.3).

This completes Case~(2).

This completes the proof.
\end{proof}

With the above propositions we can introduce the main result of this chapter.

\begin{restatable}{theorem}{TDPthFree}
\label{TDPthFree} 
The many-sorted partial $\Sigma^{\boldsymbol{\mathcal{A}}^{(2)}}$-algebras $\llbracket \mathbf{Pth}_{\boldsymbol{\mathcal{A}}^{(2)}}\rrbracket$ and $\mathbf{T}_{\boldsymbol{\mathcal{E}}^{\boldsymbol{\mathcal{A}}^{(2)}}}(\mathbf{Pth}_{\boldsymbol{\mathcal{A}}^{(2)}})$ are isomorphic.
\end{restatable}
\begin{proof}
To follow the proof, the reader is urged to consider the diagram given in Figure~\ref{FDBigPth}. By Corollary~\ref{CDVarPr}, the $\Sigma^{\boldsymbol{\mathcal{A}}^{(2)}}$-homomorphism $\mathrm{pr}^{\llbracket\cdot\rrbracket}$ from $\mathbf{Pth}_{\boldsymbol{\mathcal{A}}^{(2)}}$ to $\llbracket \mathbf{Pth}_{\boldsymbol{\mathcal{A}}^{(2)}}\rrbracket$ extends to a $\Sigma^{\boldsymbol{\mathcal{A}}^{(2)}}$-homomorphism, $\mathrm{pr}^{\llbracket\cdot\rrbracket\mathsf{p}}$ from $\mathbf{T}_{\boldsymbol{\mathcal{E}}^{\boldsymbol{\mathcal{A}}^{(2)}}}(\mathbf{Pth}_{\boldsymbol{\mathcal{A}}^{(2)}})$ to $\llbracket \mathbf{Pth}_{\boldsymbol{\mathcal{A}}^{(2)}}\rrbracket$.  By Remark~\ref{RDQPUniv}, the mapping $
(\mathrm{pr}^{\equiv^{\llbracket 2\rrbracket}}
\circ
\mathrm{ip}^{(2,X)@}
\circ
\mathrm{CH}^{(2)})^{\natural}
$ from $
\llbracket 
\mathbf{Pth}_{\boldsymbol{\mathcal{A}}^{(2)}}
\rrbracket
$ to $
\mathbf{T}_{\boldsymbol{\mathcal{E}}^{\boldsymbol{\mathcal{A}}^{(2)}}}
(
\mathbf{Pth}_{\boldsymbol{\mathcal{A}}^{(2)}}
)
$
is a $\Sigma^{\boldsymbol{\mathcal{A}}^{(2)}}$-homomorphism. 

Thus, the composition
$$
(\mathrm{pr}^{\equiv^{\llbracket 2\rrbracket}}
\circ
\mathrm{ip}^{(2,X)@}
\circ
\mathrm{CH}^{(2)})^{\natural}
\circ
\mathrm{pr}^{\llbracket\cdot\rrbracket\mathsf{p}}
$$
is a $\Sigma^{\boldsymbol{\mathcal{A}}^{(2)}}$-endomorphism of $\mathbf{T}_{\boldsymbol{\mathcal{E}}^{\boldsymbol{\mathcal{A}}^{(2)}}}(\mathbf{Pth}_{\boldsymbol{\mathcal{A}}^{(2)}})$.

Moreover, this composition satisfies the following chain of equalities
\begin{flushleft}
$(\mathrm{pr}^{\equiv^{\llbracket 2\rrbracket}}
\circ
\mathrm{ip}^{(2,X)@}
\circ
\mathrm{CH}^{(2)})^{\natural}
\circ
\mathrm{pr}^{\llbracket\cdot\rrbracket\mathsf{p}}
\circ
\eta^{(\llbracket 2\rrbracket,\mathbf{Pth}_{\boldsymbol{\mathcal{A}}^{(2)}})}$
\allowdisplaybreaks
\begin{align*}
\qquad
&=
(\mathrm{pr}^{\equiv^{\llbracket 2\rrbracket}}
\circ
\mathrm{ip}^{(2,X)@}
\circ
\mathrm{CH}^{(2)})^{\natural}
\circ
\mathrm{pr}^{\llbracket\cdot\rrbracket\mathsf{p}}
\tag{1}
\\
&=
\mathrm{pr}^{\equiv^{\llbracket 2\rrbracket}}
\circ
\mathrm{ip}^{(2,X)@}
\circ
\mathrm{CH}^{(2)}
\tag{2}
\\
&=
\eta^{(\llbracket 2\rrbracket,\mathbf{Pth}_{\boldsymbol{\mathcal{A}}^{(2)}})}.
\tag{3}
\end{align*}
\end{flushleft}

In the just stated chain of equalities, the first equality holds by Corollary~\ref{CDVarPr}; the second equality follows from Remark~\ref{RDQPUniv}; finally the last equality holds by Proposition~\ref{PDQEta}.

Let us recall that the identity mapping $\mathrm{id}^{\mathbf{T}_{\boldsymbol{\mathcal{E}}^{\boldsymbol{\mathcal{A}}^{(2)}}}(\mathbf{Pth}_{\boldsymbol{\mathcal{A}}^{(2)}})}$ is also a $\Sigma^{\boldsymbol{\mathcal{A}}^{(2)}}$-endomorphism of $\mathbf{T}_{\boldsymbol{\mathcal{E}}^{\boldsymbol{\mathcal{A}}^{(2)}}}(\mathbf{Pth}_{\boldsymbol{\mathcal{A}}^{(2)}})$ satisfying that
$$
\mathrm{id}^{\mathbf{T}_{\boldsymbol{\mathcal{E}}^{\boldsymbol{\mathcal{A}}^{(2)}}}(\mathbf{Pth}_{\boldsymbol{\mathcal{A}}^{(2)}})}
\circ
\eta^{(\llbracket 2\rrbracket,\mathbf{Pth}_{\boldsymbol{\mathcal{A}}^{(2)}})}
=\eta^{(\llbracket 2\rrbracket,\mathbf{Pth}_{\boldsymbol{\mathcal{A}}^{(2)}})}.
$$

Hence, by the universal property of $\mathbf{T}_{\boldsymbol{\mathcal{E}}^{\boldsymbol{\mathcal{A}}^{(2)}}}(\mathbf{Pth}_{\boldsymbol{\mathcal{A}}^{(2)}})$, we have 
$$
(\mathrm{pr}^{\equiv^{\llbracket 2\rrbracket}}
\circ
\mathrm{ip}^{(2,X)@}
\circ
\mathrm{CH}^{(2)})^{\natural}
\circ
\mathrm{pr}^{\llbracket\cdot\rrbracket\mathsf{p}}
=
\mathrm{id}^{\mathbf{T}_{\boldsymbol{\mathcal{E}}^{\boldsymbol{\mathcal{A}}^{(2)}}}(\mathbf{Pth}_{\boldsymbol{\mathcal{A}}^{(2)}})}.
$$

Therefore, $(\mathrm{pr}^{\equiv^{\llbracket 2\rrbracket}}
\circ
\mathrm{ip}^{(2,X)@}
\circ
\mathrm{CH}^{(2)})^{\natural}$ is a $\Sigma^{\boldsymbol{\mathcal{A}}^{(2)}}$-isomorphism from  $\llbracket \mathbf{Pth}_{\boldsymbol{\mathcal{A}}^{(2)}}\rrbracket$ to $\mathbf{T}_{\boldsymbol{\mathcal{E}}^{\boldsymbol{\mathcal{A}}^{(2)}}}(\mathbf{Pth}_{\boldsymbol{\mathcal{A}}^{(2)}})$ whose inverse is $\mathrm{pr}^{\llbracket\cdot\rrbracket\mathsf{p}}$.

This finishes the proof.
\end{proof}

\begin{restatable}{corollary}{CDPTFree}
\label{CDPTFree} 
$\llbracket\mathbf{PT}_{\boldsymbol{\mathcal{A}}^{(2)}}\rrbracket$ is a many-sorted partial  $\Sigma^{\boldsymbol{\mathcal{A}}^{(2)}}$-algebra in $\mathbf{PAlg}(\boldsymbol{\mathcal{E}}^{\boldsymbol{\mathcal{A}}^{(2)}})$. Indeed, the many-sorted partial $\Sigma^{\boldsymbol{\mathcal{A}}^{(2)}}$-algebras $\llbracket\mathbf{PT}_{\boldsymbol{\mathcal{A}}^{(2)}}\rrbracket$ and $\mathbf{T}_{\boldsymbol{\mathcal{E}}^{\boldsymbol{\mathcal{A}}^{(2)}}}(\mathbf{Pth}_{\boldsymbol{\mathcal{A}}^{(2)}})$ are isomorphic.
\end{restatable}

\begin{figure}
\begin{center}
\begin{tikzpicture}
[ACliment/.style={-{To [angle'=45, length=5.75pt, width=4pt, round]}
}, scale=0.6]
\node[] (P3) 		at 	(0,-4) 	[] 	{$\mathbf{Pth}_{\boldsymbol{\mathcal{A}}^{(2)}}$};
\node[] (TEP) 	at 	(6,0) 	[] 	{$\mathbf{T}_{\boldsymbol{\mathcal{E}}^{\boldsymbol{\mathcal{A}}^{(2)}}}(\mathbf{Pth}_{\boldsymbol{\mathcal{A}}^{(2)}})$};
\node[] (PK) 		at 	(6,-4) 	[] 	{$\llbracket \mathbf{Pth}_{\boldsymbol{\mathcal{A}}^{(2)}}\rrbracket$};
\node[] (TEP2) 	at 	(6,-8) 	[] 	{$\mathbf{T}_{\boldsymbol{\mathcal{E}}^{\boldsymbol{\mathcal{A}}^{(2)}}}
(\mathbf{Pth}_{\boldsymbol{\mathcal{A}}^{(2)}})$};
\draw[ACliment, bend left]  (P3) 	to node [above left]	
{$\eta^{(\llbracket 2 \rrbracket,\mathbf{Pth}_{\boldsymbol{\mathcal{A}}^{(2)}})}$} (TEP);
\draw[ACliment, bend right]  (P3) 	to node [below left]	
{$\eta^{(\llbracket 2 \rrbracket,\mathbf{Pth}_{\boldsymbol{\mathcal{A}}^{(2)}})}$} (TEP2);
\draw[ACliment]  (P3) 	to node [above]	
{$\mathrm{pr}^{\llbracket\cdot \rrbracket}$} (PK);
\draw[ACliment]  (PK) 	to node [midway, fill=white]	
{$(\mathrm{pr}^{\equiv^{\llbracket 2\rrbracket}}
\circ
\mathrm{ip}^{(2,X)@}
\circ
\mathrm{CH}^{(2)})^{\natural}$} (TEP2);
\draw[ACliment]  (TEP) 	to node [right]	
{$\mathrm{pr}^{\llbracket\cdot \rrbracket\mathsf{p}}$} (PK);
\draw[ACliment, rounded corners] (TEP.east)
--
 ($(TEP.east)+(1.5,0)$)
to node [right]	
{$\mathrm{id}^{\mathbf{T}_{\boldsymbol{\mathcal{E}}^{\boldsymbol{\mathcal{A}}^{(2)}}}(\mathbf{Pth}_{\boldsymbol{\mathcal{A}}^{(2)}})}$} 
($(TEP2.east)+(1.5,0)$)
-- (TEP2.east);
\end{tikzpicture}
\end{center}
\caption{The maps in Theorem~\ref{TDPthFree}.}
\label{FDBigPth}
\end{figure}

\chapter{Second-order translations for second-order path terms}\label{S2P}

In this chapter, the last chapter of the second part, introduces the notions of second-order  elementary translation and second-order translation. It is shown that second-order translations are well-behaved, in the sense that if two second-order path terms $P$ and $Q$ in $\mathrm{PT}_{\boldsymbol{\mathcal{A}}^{(2)},s}$ are such that $\mathrm{ip}^{(2,X)@}_{s}(P)$ and $\mathrm{ip}^{(2,X)@}_{s}(Q)$ have the same $([1],2)$-source and $([1],2)$-target, then the terms that result from applying a second-order translation to $P$ and $Q$ are also second-order path terms and, when turned into second-order paths by means of $\mathrm{ip}^{(2,X)@}$, they still have the same $([1],2)$-source and $([1],2)$-target. Thus, the same $(0,2)$-source and $(0,2)$-target. Moreover, it is shown that, for two second-order path terms $P$ and $Q$ in $\mathrm{PT}_{\boldsymbol{\mathcal{A}}^{(2)},s}$ that are $\Theta^{\llbracket 2 \rrbracket}_{s}$-related, then its respective translations are also $\Theta^{\llbracket 2 \rrbracket}_{s}$-related.


We next define for the many-sorted partial $\Sigma^{\boldsymbol{\mathcal{A}}^{(2)}}$-algebra $\mathbf{T}_{\Sigma^{\boldsymbol{\mathcal{A}}^{(2)}}}(X)$ the concepts of second-order elementary  translation and of second-order translation respect to it.

\begin{restatable}{definition}{DDETrans}
\label{DDETrans} 
\index{translation!second-order!elementary}
Let $t$ be a sort in $S$. We will denote by $\mathrm{Etl}_{t}(\mathbf{T}_{\Sigma^{\boldsymbol{\mathcal{A}}^{(2)}}}(X))$ the subset $(\mathrm{Etl}_{t}(\mathbf{T}_{\Sigma^{\boldsymbol{\mathcal{A}}^{(2)}}}(X))_{s})_{s\in S}$ of $(\mathrm{Hom}(\mathrm{T}_{\Sigma^{\boldsymbol{\mathcal{A}}^{(2)}}}(X)_{t},\mathrm{T}_{\Sigma^{\boldsymbol{\mathcal{A}}^{(2)}}}(X)_{s}))_{s\in S}$ defined, for every $s\in S$, as follows: for every mapping $T^{(2)}\in \mathrm{Hom}(\mathrm{T}_{\Sigma^{\boldsymbol{\mathcal{A}}^{(2)}}}(X)_{t},\mathrm{T}_{\Sigma^{\boldsymbol{\mathcal{A}}^{(2)}}}(X)_{s})$, $T^{(2)}\in \mathrm{Etl}_{t}(\mathbf{T}_{\Sigma^{\boldsymbol{\mathcal{A}}^{(2)}}}(X))_{s}$ if and only if one of the following conditions holds
\begin{enumerate}
\item There is a word $\mathbf{s}\in S^{\star}-\{\lambda\}$, an index $k\in\bb{\mathbf{s}}$, an operation symbol $\sigma\in \Sigma_{\mathbf{s},s}$, a family of second-order paths $(\mathfrak{P}^{(2)}_{j})_{j\in k}\in \prod_{j\in k}\mathrm{Pth}_{\boldsymbol{\mathcal{A}}^{(2)},s_{j}}$ and a family of second-order paths $(\mathfrak{P}^{(2)}_{l})_{l\in \bb{\mathbf{s}}-(k+1)}\in \prod_{l\in \bb{\mathbf{s}}-(k+1)}\mathrm{Pth}_{\boldsymbol{\mathcal{A}}^{(2)},s_{l}}$ (recall that $k+1=\{0,\cdots, k\}$ and that $\bb{\mathbf{s}}-(k+1)=\{k+1,\cdots, \bb{\mathbf{s}}-1\}$) such that $s_{k}=t$ and, for every $P\in\mathrm{T}_{\Sigma^{\boldsymbol{\mathcal{A}}^{(2)}}}(X)_{t}$
\begin{multline*}
T^{(2)}(P)=\sigma^{\mathbf{T}_{\Sigma^{\boldsymbol{\mathcal{A}}^{(2)}}}(X)}\left(
\mathrm{CH}^{(2)}_{s_{0}}\left(\mathfrak{P}^{(2)}_{0}\right),
\cdots,
\mathrm{CH}^{(2)}_{s_{k-1}}\left(\mathfrak{P}^{(2)}_{k-1}\right),
\right.
\\
\left.
P,
\mathrm{CH}^{(2)}_{s_{k+1}}\left(\mathfrak{P}^{(2)}_{k+1}\right),
\cdots,
\mathrm{CH}^{(2)}_{s_{\bb{\mathbf{s}}-1}}\left(\mathfrak{P}^{(2)}_{\bb{\mathbf{s}}-1}\right)
\right);
\end{multline*}
\item It holds that $s=t$ and there is a second-order path $\mathfrak{P}^{(2)}\in\mathrm{Pth}_{\boldsymbol{\mathcal{A}}^{(2)},s}$ such that, for every $P\in\mathrm{T}_{\Sigma^{\boldsymbol{\mathcal{A}}^{(2)}}}(X)_{s}$,
\begin{align*}
T^{(2)}(P)&=P\circ^{0\mathbf{T}_{\Sigma^{\boldsymbol{\mathcal{A}}^{(2)}}}(X)}_{s}\mathrm{CH}^{(2)}_{s}\left(\mathfrak{P}^{(2)}\right); \mbox{ or }
\\
T^{(2)}(P)&=\mathrm{CH}^{(2)}_{s}\left(\mathfrak{P}^{(2)}\right)\circ^{0\mathbf{T}_{\Sigma^{\boldsymbol{\mathcal{A}}^{(2)}}}(X)}_{s}P.
\end{align*}
\item It holds that $s=t$ and there is a second-order path $\mathfrak{P}^{(2)}\in\mathrm{Pth}_{\boldsymbol{\mathcal{A}}^{(2)},s}$ such that, for every $P\in\mathrm{T}_{\Sigma^{\boldsymbol{\mathcal{A}}^{(2)}}}(X)_{s}$,
\begin{align*}
T^{(2)}(P)&=P\circ^{1\mathbf{T}_{\Sigma^{\boldsymbol{\mathcal{A}}^{(2)}}}(X)}_{s}\mathrm{CH}^{(2)}_{s}\left(\mathfrak{P}^{(2)}\right); \mbox{ or }
\\
T^{(2)}(P)&=\mathrm{CH}^{(2)}_{s}\left(\mathfrak{P}^{(2)}\right)\circ^{1\mathbf{T}_{\Sigma^{\boldsymbol{\mathcal{A}}^{(2)}}}(X)}_{s}P.
\end{align*}
\end{enumerate}

We will sometimes use the following presentations of second-order elementary translations, adding an underlined space to denote where the variable will be placed.
\begin{multline*}
T^{(2)}=\sigma^{\mathbf{T}_{\Sigma^{\boldsymbol{\mathcal{A}}^{(2)}}}(X)}\left(
\mathrm{CH}^{(2)}_{s_{0}}\left(\mathfrak{P}^{(2)}_{0}\right),
\cdots,
\mathrm{CH}^{(2)}_{s_{k-1}}\left(\mathfrak{P}^{(2)}_{k-1}\right),
\right.
\\
\left.
\underline{\quad},
\mathrm{CH}^{(2)}_{s_{k+1}}\left(\mathfrak{P}^{(2)}_{k+1}\right),
\cdots,
\mathrm{CH}^{(2)}_{s_{\bb{\mathbf{s}}-1}}\left(\mathfrak{P}^{(2)}_{\bb{\mathbf{s}}-1}\right)
\right);
\tag{1}
\end{multline*}
\begin{align*}
T^{(2)}&=\underline{\quad}\circ^{0\mathbf{T}_{\Sigma^{\boldsymbol{\mathcal{A}}^{(2)}}}(X)}_{s}\mathrm{CH}^{(2)}_{s}\left(\mathfrak{P}^{(2)}\right); \mbox{ or }
\\
T^{(2)}&=\mathrm{CH}^{(2)}_{s}\left(\mathfrak{P}^{(2)}\right)\circ^{0\mathbf{T}_{\Sigma^{\boldsymbol{\mathcal{A}}^{(2)}}}(X)}_{s}\underline{\quad};
\tag{2}
\end{align*}
\begin{align*}
T^{(2)}&=\underline{\quad}\circ^{1\mathbf{T}_{\Sigma^{\boldsymbol{\mathcal{A}}^{(2)}}}(X)}_{s}\mathrm{CH}^{(2)}_{s}\left(\mathfrak{P}^{(2)}\right); \mbox{ or }
\\
T^{(2)}&=\mathrm{CH}^{(2)}_{s}\left(\mathfrak{P}^{(2)}\right)\circ^{1\mathbf{T}_{\Sigma^{\boldsymbol{\mathcal{A}}^{(2)}}}(X)}_{s}\underline{\quad}.
\tag{3}
\end{align*}

In the first case we will say that $T^{(2)}$ is a \emph{second-order elementary translation of type $\sigma$}, in the second case we will say that $T^{(2)}$ is an \emph{second-order elementary translation of type $\circ^{0}_{s}$} and, finally, in the third case we will say that $T^{(2)}$ is an \emph{second-order elementary translation of type $\circ^{1}_{s}$}. We will call the elements of $\mathrm{Etl}_{t}(\mathbf{T}_{\Sigma^{\boldsymbol{\mathcal{A}}^{(2)}}}(X))_{s}$ the \emph{second-order $t$-elementary translations of sort $s$ for $\mathbf{T}_{\Sigma^{\boldsymbol{\mathcal{A}}^{(2)}}}(X)$}.
\end{restatable}

\begin{restatable}{definition}{DDTrans}
\label{DDTrans} 
\index{translation!second-order}
Let $t$ be a sort in $S$. We will denote by $\mathrm{Tl}_{t}(\mathbf{T}_{\Sigma^{\boldsymbol{\mathcal{A}}^{(2)}}}(X))$ the subset $(\mathrm{Tl}_{t}(\mathbf{T}_{\Sigma^{\boldsymbol{\mathcal{A}}^{(2)}}}(X))_{s})_{s\in S}$ of $(\mathrm{Hom}(\mathrm{T}_{\Sigma^{\boldsymbol{\mathcal{A}}^{(2)}}}(X)_{t},\mathrm{T}_{\Sigma^{\boldsymbol{\mathcal{A}}^{(2)}}}(X)_{s}))_{s\in S}$ defined, for every $s\in S$, as follows: for every mapping $T^{(2)}\in \mathrm{Hom}(\mathrm{T}_{\Sigma^{\boldsymbol{\mathcal{A}}^{(2)}}}(X)_{t},\mathrm{T}_{\Sigma^{\boldsymbol{\mathcal{A}}^{(2)}}}(X)_{s})$, $T^{(2)}\in \mathrm{Tl}_{t}(\mathbf{T}_{\Sigma^{\boldsymbol{\mathcal{A}}^{(2)}}}(X))_{s}$ if and only if there is an $n\in \mathbb{N}-1$, a word $\mathbf{w}\in S^{n+1}$, and a family $(T^{(2)}_{j})_{j\in n}$ such that $w_{0}=t$, $w_{n}=s$, and, for every $j\in n$, $T^{(2)}_{j}\in \mathrm{Etl}_{w_{j}}(\mathbf{T}_{\Sigma^{\boldsymbol{\mathcal{A}}^{(2)}}}(X))_{w_{j+1}}$ and 
$T^{(2)}=T^{(2)}_{n-1}\circ \cdots \circ T^{(2)}_{0}$. We will sometimes refer to the composition $T^{(2)}_{n-2}\circ \cdots \circ T^{(2)}_{0}$ as the \emph{derived second-order translation} of $T^{(2)}$, denoted by $T^{(2)'}$, and we will express $T^{(2)}$ as $T^{(2)}_{n-1}\circ T^{(2)'}$ or in  operator form
\begin{multline*}
T^{(2)}=\sigma^{\mathbf{T}_{\Sigma^{\boldsymbol{\mathcal{A}}^{(2)}}}(X)}\left(
\mathrm{CH}^{(2)}_{s_{0}}\left(\mathfrak{P}^{(2)}_{0}\right),
\cdots,
\mathrm{CH}^{(2)}_{s_{k-1}}\left(\mathfrak{P}^{(2)}_{k-1}\right),
\right.
\\
\left.
T^{(2)'},
\mathrm{CH}^{(2)}_{s_{k+1}}\left(\mathfrak{P}^{(2)}_{k+1}\right),
\cdots,
\mathrm{CH}^{(2)}_{s_{\bb{\mathbf{s}}-1}}\left(\mathfrak{P}^{(2)}_{\bb{\mathbf{s}}-1}\right)
\right);
\tag{1}
\end{multline*}
\begin{align*}
T^{(2)}&=T^{(2)'}\circ^{0\mathbf{T}_{\Sigma^{\boldsymbol{\mathcal{A}}^{(2)}}}(X)}_{s}\mathrm{CH}^{(2)}_{s}\left(\mathfrak{P}^{(2)}\right);
\\
T^{(2)}&=\mathrm{CH}^{(2)}_{s}\left(\mathfrak{P}^{(2)}\right)\circ^{0\mathbf{T}_{\Sigma^{\boldsymbol{\mathcal{A}}^{(2)}}}(X)}_{s}T^{(2)'};
\tag{2}
\end{align*}
\begin{align*}
T^{(2)}&=T^{(2)'}\circ^{1\mathbf{T}_{\Sigma^{\boldsymbol{\mathcal{A}}^{(2)}}}(X)}_{s}\mathrm{CH}^{(2)}_{s}\left(\mathfrak{P}^{(2)}\right);
\\
T^{(2)}&=\mathrm{CH}^{(2)}_{s}\left(\mathfrak{P}^{(2)}\right)\circ^{1\mathbf{T}_{\Sigma^{\boldsymbol{\mathcal{A}}^{(2)}}}(X)}_{s}T^{(2)'}.
\tag{3}
\end{align*}

This operator presentation means that, for every $P\in\mathrm{T}_{\Sigma^{\boldsymbol{\mathcal{A}}^{(2)}}}(X)_{t}$, $T^{(2)}(P)=T^{(2)}_{n-1}(T^{(2)'}(P))$. The underlined space notation can be extended to second-order translations as well. We will say that $T^{(2)}$ is a \emph{second-order translation of type $\sigma$}, a \emph{second-order translation of type $\circ^{0}_{s}$}, or a \emph{second-order translation of type $\circ^{1}_{s}$}, if $\mathrm{T}^{(2)}_{n-1}$ is a second-order elementary translation of type $\sigma$, a second-order elementary translation of type $\circ^{0}_{s}$, or a second-order elementary translation of type $\circ^{1}_{s}$, respectively. 

We will call $n$ the \emph{height of $T^{(2)}$} and we will denote this fact by $\bb{T^{(2)}}=n$. In this regard, second-order elementary translation have height $1$ and, if $T^{(2)}$ is a second-order translation of height $n$, i.e., $\bb{T^{(2)}}=n$, then its derived second-order translation has height $n-1$, i.e., $\bb{T^{(2)'}}=n-1$.

We will call the elements of $\mathrm{Tl}_{t}(\mathbf{T}_{\Sigma^{\boldsymbol{\mathcal{A}}^{(2)}}}(X))_{s}$ the \emph{second-order $t$-translations of sort $s$ for $\mathbf{T}_{\Sigma^{\boldsymbol{\mathcal{A}}^{(2)}}}(X)$}. Besides, for every $t\in S$, the mapping $\mathrm{id}^{\mathrm{T}_{\Sigma^{\boldsymbol{\mathcal{A}}^{(2)}}}(X)_{t}}$ will be viewed as an element of $\mathrm{Tl}_{t}(\mathbf{T}_{\Sigma^{\boldsymbol{\mathcal{A}}^{(2)}}}(X))_{t}$. The identity second-order translation has no associated type and we will consider that it has height $0$, i.e., $\bb{\mathrm{id}^{\mathrm{T}_{\Sigma^{\boldsymbol{\mathcal{A}}^{(2)}}}(X)_{t}}}=0$.
Viewing the identity as a second-order translation allows us to give an operator presentation for second-order elementary translations, where the derived second-order translation is the identity.
\end{restatable}

\begin{restatable}{lemma}{LDTransWD}
\label{LDTransWD} 
Let $s, t$ be sorts in $S$, $T^{(2)}$ a second-order translation in 
$\mathrm{Tl}_{t}(\mathrm{T}_{\Sigma^{\boldsymbol{\mathcal{A}}^{(2)}}}(X))_{s}$ and let $M$ and $N$ be second-order path terms in
$\mathrm{PT}_{\boldsymbol{\mathcal{A}}^{(2)},t}$ satisfying that 
 $$
\left(
\mathrm{ip}^{(2,X)@}_{t}\left(
M
\right), 
\mathrm{ip}^{(2,X)@}_{t}\left(
N
\right)\right)
\in\mathrm{Ker}\left(
\mathrm{sc}^{([1],2)}
\right)_{t}
\cap\mathrm{Ker}\left(
\mathrm{tg}^{([1],2)}
\right)_{t}.
$$
Then the following properties hold 
\begin{itemize}
\item[(i)]
$T^{(2)}(M)$ is a second-order path term in $\mathrm{PT}_{\boldsymbol{\mathcal{A}}^{(2)},s}$ if, and only if, $T^{(2)}(N)$ is a second-order path term in  $\mathrm{PT}_{\boldsymbol{\mathcal{A}}^{(2)},s}$;
\item[(ii)] If either $T^{(2)}(M)$ or $T^{(2)}(N)$ is a second-order path term in $\mathrm{PT}_{\boldsymbol{\mathcal{A}}^{(2)},s}$, then 
\[
\left(
\mathrm{ip}^{(2,X)@}_{s}\left(
T^{(2)}(M)\right),
\mathrm{ip}^{(2,X)@}_{s}\left(
T^{(2)}(N)\right)
\right)
\in\mathrm{Ker}\left(
\mathrm{sc}^{([1],2)}
\right)_{s}
\cap\mathrm{Ker}\left(
\mathrm{tg}^{([1],2)}
\right)_{s}.
\]
\end{itemize}
\end{restatable}
\begin{proof}
We will prove the statement by induction on the height of the second-order translation. 

\textsf{Base case of the induction.}

The statement follows easily when $s=t$ and $T^{(2)}=\mathrm{id}^{\mathrm{T}_{\Sigma^{\boldsymbol{\mathcal{A}}^{(2)}}}(X)_{s}}$.

We consider the different possibilities for the case of $T^{(2)}$ being a second-order elementary translation according to Definition~\ref{DDETrans}.

Assume that there is a word $\mathbf{s}\in S^{\star}-\{\lambda\}$, an index $k\in\bb{\mathbf{s}}$, an operation symbol $\sigma\in \Sigma_{\mathbf{s},s}$, a family of second-order paths $(\mathfrak{P}^{(2)}_{j})_{j\in k}\in \prod_{j\in k}\mathrm{Pth}_{\boldsymbol{\mathcal{A}}^{(2)},s_{j}}$ and a family of paths $(\mathfrak{P}_{l})_{l\in \bb{\mathbf{s}}-(k+1)}\in \prod_{l\in \bb{\mathbf{s}}-(k+1)}\mathrm{Pth}_{\boldsymbol{\mathcal{A}}^{(2)},s_{l}}$ such that $s_{k}=t$ and
\begin{multline*}
T^{(2)}=\sigma^{\mathbf{T}_{\Sigma^{\boldsymbol{\mathcal{A}}^{(2)}}}(X)}\left(
\mathrm{CH}^{(2)}_{s_{0}}\left(\mathfrak{P}^{(2)}_{0}\right),
\cdots,
\mathrm{CH}^{(2)}_{s_{k-1}}\left(\mathfrak{P}^{(2)}_{k-1}\right),
\right.
\\
\left.
\underline{\quad},
\mathrm{CH}^{(2)}_{s_{k+1}}\left(\mathfrak{P}^{(2)}_{k+1}\right),
\cdots,
\mathrm{CH}^{(2)}_{s_{\bb{\mathbf{s}}-1}}\left(\mathfrak{P}^{(2)}_{\bb{\mathbf{s}}-1}\right)
\right).
\end{multline*}

Note that, for this case, the terms $T^{(2)}(M)$ and $T^{(2)}(N)$ are always second-order path terms. Therefore, the first condition of the proposition holds.

Now, regarding the $([1],2)$-source, the following chain of equalities holds
\begin{flushleft}
$\mathrm{sc}^{([1],2)}_{s}\left(\mathrm{ip}^{(2,X)@}_{s}\left(
T^{(2)}(M)\right)\right)$
\allowdisplaybreaks
\begin{align*}
\qquad&=
\mathrm{sc}^{([1],2)}_{s}\left(
\mathrm{ip}^{(2,X)@}_{s}\left(
\sigma^{\mathbf{T}_{\Sigma^{\boldsymbol{\mathcal{A}}^{(2)}}}(X)}\left(
\mathrm{CH}^{(2)}_{s_{0}}\left(\mathfrak{P}^{(2)}_{0}\right),
\cdots,
\mathrm{CH}^{(2)}_{s_{k-1}}\left(\mathfrak{P}^{(2)}_{k-1}\right),
\right.\right.\right.
\\&\qquad\qquad\qquad\qquad\qquad
\left.\left.\left.
M,
\mathrm{CH}^{(2)}_{s_{k+1}}\left(\mathfrak{P}^{(2)}_{k+1}\right),
\cdots,
\mathrm{CH}^{(2)}_{s_{\bb{\mathbf{s}}-1}}\left(\mathfrak{P}^{(2)}_{\bb{\mathbf{s}}-1}\right)
\right)
\right)\right)
\tag{1}
\\&=
\mathrm{sc}^{([1],2)}_{s}\left(
\sigma^{\mathbf{Pth}_{\boldsymbol{\mathcal{A}}^{(2)}}}\left(
\mathrm{ip}^{(2,X)@}_{s_{0}}\left(
\mathrm{CH}^{(2)}_{s_{0}}\left(\mathfrak{P}^{(2)}_{0}\right)\right),
\cdots,
\right.\right.
\\&\qquad\qquad
\mathrm{ip}^{(2,X)@}_{s_{k-1}}\left(
\mathrm{CH}^{(2)}_{s_{k-1}}\left(\mathfrak{P}^{(2)}_{k-1}\right)\right),
\\&\qquad\qquad\qquad\qquad
\mathrm{ip}^{(2,X)@}_{t}\left(
M\right),
\mathrm{ip}^{(2,X)@}_{s_{k+1}}\left(
\mathrm{CH}^{(2)}_{s_{k+1}}\left(\mathfrak{P}^{(2)}_{k+1}\right)\right),
\\&\qquad\qquad\qquad\qquad\qquad\qquad\qquad\quad
\left.\left.
\cdots,
\mathrm{ip}^{(2,X)@}_{s_{\bb{\mathbf{s}}-1}}\left(
\mathrm{CH}^{(2)}_{s_{\bb{\mathbf{s}}-1}}\left(\mathfrak{P}^{(2)}_{\bb{\mathbf{s}}-1}\right)\right)
\right)
\right)
\tag{2}
\\&=
\sigma^{[\mathbf{PT}_{\boldsymbol{\mathcal{A}}}]}\left(
\mathrm{sc}^{([1],2)}_{s_{0}}\left(
\mathrm{ip}^{(2,X)@}_{s_{0}}\left(
\mathrm{CH}^{(2)}_{s_{0}}\left(\mathfrak{P}^{(2)}_{0}\right)\right)\right),
\cdots,
\right.
\\&\qquad\qquad
\mathrm{sc}^{([1],2)}_{s_{k-1}}\left(
\mathrm{ip}^{(2,X)@}_{s_{k-1}}\left(
\mathrm{CH}^{(2)}_{s_{k-1}}\left(\mathfrak{P}^{(2)}_{k-1}\right)\right)\right),
\\&\qquad\qquad\qquad
\mathrm{sc}^{([1],2)}_{t}\left(
\mathrm{ip}^{(2,X)@}_{t}\left(
M\right)\right),
\mathrm{sc}^{([1],2)}_{s_{k+1}}\left(
\mathrm{ip}^{(2,X)@}_{s_{k+1}}\left(
\mathrm{CH}^{(2)}_{s_{k+1}}\left(\mathfrak{P}^{(2)}_{k+1}\right)\right)\right),
\\&\qquad\qquad\qquad\qquad\qquad\qquad
\left.
\cdots,
\mathrm{sc}^{([1],2)}_{s_{\bb{\mathbf{s}}-1}}\left(
\mathrm{ip}^{(2,X)@}_{s_{\bb{\mathbf{s}}-1}}\left(
\mathrm{CH}^{(2)}_{s_{\bb{\mathbf{s}}-1}}\left(\mathfrak{P}^{(2)}_{\bb{\mathbf{s}}-1}\right)\right)
\right)
\right)
\tag{3}
\\&=
\sigma^{[\mathbf{PT}_{\boldsymbol{\mathcal{A}}}]}\left(
\mathrm{sc}^{([1],2)}_{s_{0}}\left(
\mathrm{ip}^{(2,X)@}_{s_{0}}\left(
\mathrm{CH}^{(2)}_{s_{0}}\left(\mathfrak{P}^{(2)}_{0}\right)\right)\right),
\cdots,
\right.
\\&\qquad\qquad
\mathrm{sc}^{([1],2)}_{s_{k-1}}\left(
\mathrm{ip}^{(2,X)@}_{s_{k-1}}\left(
\mathrm{CH}^{(2)}_{s_{k-1}}\left(\mathfrak{P}^{(2)}_{k-1}\right)\right)\right),
\\&\qquad\qquad\qquad
\mathrm{sc}^{([1],2)}_{t}\left(
\mathrm{ip}^{(2,X)@}_{t}\left(
N\right)\right),
\mathrm{sc}^{([1],2)}_{s_{k+1}}\left(
\mathrm{ip}^{(2,X)@}_{s_{k+1}}\left(
\mathrm{CH}^{(2)}_{s_{k+1}}\left(\mathfrak{P}^{(2)}_{k+1}\right)\right)\right),
\\&\qquad\qquad\qquad\qquad\qquad\qquad
\left.
\cdots,
\mathrm{sc}^{([1],2)}_{s_{\bb{\mathbf{s}}-1}}\left(
\mathrm{ip}^{(2,X)@}_{s_{\bb{\mathbf{s}}-1}}\left(
\mathrm{CH}^{(2)}_{s_{\bb{\mathbf{s}}-1}}\left(\mathfrak{P}^{(2)}_{\bb{\mathbf{s}}-1}\right)\right)
\right)
\right)
\tag{4}
\\&=
\mathrm{sc}^{([1],2)}_{s}\left(
\sigma^{\mathbf{Pth}_{\boldsymbol{\mathcal{A}}^{(2)}}}\left(
\mathrm{ip}^{(2,X)@}_{s_{0}}\left(
\mathrm{CH}^{(2)}_{s_{0}}\left(\mathfrak{P}^{(2)}_{0}\right)\right),
\cdots,
\right.\right.
\\&\qquad\qquad
\mathrm{ip}^{(2,X)@}_{s_{k-1}}\left(
\mathrm{CH}^{(2)}_{s_{k-1}}\left(\mathfrak{P}^{(2)}_{k-1}\right)\right),
\\&\qquad\qquad\qquad\qquad
\mathrm{ip}^{(2,X)@}_{t}\left(
N\right),
\mathrm{ip}^{(2,X)@}_{s_{k+1}}\left(
\mathrm{CH}^{(2)}_{s_{k+1}}\left(\mathfrak{P}^{(2)}_{k+1}\right)\right),
\\&\qquad\qquad\qquad\qquad\qquad\qquad\qquad\quad
\left.\left.
\cdots,
\mathrm{ip}^{(2,X)@}_{s_{\bb{\mathbf{s}}-1}}\left(
\mathrm{CH}^{(2)}_{s_{\bb{\mathbf{s}}-1}}\left(\mathfrak{P}^{(2)}_{\bb{\mathbf{s}}-1}\right)\right)
\right)
\right)
\tag{5}
\\&
\mathrm{sc}^{([1],2)}_{s}\left(
\mathrm{ip}^{(2,X)@}_{s}\left(
\sigma^{\mathbf{T}_{\Sigma^{\boldsymbol{\mathcal{A}}^{(2)}}}(X)}\left(
\mathrm{CH}^{(2)}_{s_{0}}\left(\mathfrak{P}^{(2)}_{0}\right),
\cdots,
\mathrm{CH}^{(2)}_{s_{k-1}}\left(\mathfrak{P}^{(2)}_{k-1}\right),
\right.\right.\right.
\\&\qquad\qquad\qquad\qquad\qquad
\left.\left.\left.
N,
\mathrm{CH}^{(2)}_{s_{k+1}}\left(\mathfrak{P}^{(2)}_{k+1}\right),
\cdots,
\mathrm{CH}^{(2)}_{s_{\bb{\mathbf{s}}-1}}\left(\mathfrak{P}^{(2)}_{\bb{\mathbf{s}}-1}\right)
\right)
\right)\right)
\tag{6}
\\&=\mathrm{sc}^{([1],2)}_{s}\left(\mathrm{ip}^{(2,X)@}_{s}\left(
T^{(2)}(N)\right)\right).
\tag{7}
\end{align*}
\end{flushleft}

In the just stated chain of equalities, the first equality unravels the description of the second-order elementary translation $T^{(2)}$; the second equality follows from the fact that $\mathrm{ip}^{(2,X)@}$ is a $\Sigma^{\boldsymbol{\mathcal{A}}^{(2)}}$-homomorphism, according to Definition~\ref{DDIp}; the third equality follows from the fact that $\mathrm{sc}^{([1],2)}$ is a $\Sigma^{\boldsymbol{\mathcal{A}}}$-homomorphism, according to Proposition~\ref{PDUCatHom}; the fourth equality follows from the fact that, by hypothesis, $
\mathrm{sc}^{([1],2)}_{t}(
\mathrm{ip}^{(2,X)@}_{t}(
M))
=
\mathrm{sc}^{([1],2)}_{t}(
\mathrm{ip}^{(2,X)@}_{t}(
N))$;
the fifth equality follows from the fact that $\mathrm{sc}^{([1],2)}$ is a $\Sigma^{\boldsymbol{\mathcal{A}}}$-homomorphism, according to Proposition~\ref{PDUCatHom}; the sixth equality follows from the fact that $\mathrm{ip}^{(2,X)@}$ is a $\Sigma^{\boldsymbol{\mathcal{A}}^{(2)}}$-homomorphism, according to Definition~\ref{DDIp}; finally, the last equality recovers  the description of the second-order elementary translation $T^{(2)}$.

The case of the $([1],2)$-target follows by a similar argument.

Now, regarding the other options for the second-order elementary translation $T^{(2)}$, it could be the case that $s=t$ and there is a second-order path $\mathfrak{P}^{(2)}\in\mathrm{Pth}_{\boldsymbol{\mathcal{A}}^{(2)},s}$ such that
\[
T^{(2)}=\underline{\quad}\circ^{0\mathbf{T}_{\Sigma^{\boldsymbol{\mathcal{A}}^{(2)}}}(X)}_{s}\mathrm{CH}^{(2)}_{s}\left(\mathfrak{P}^{(2)}\right).
\]

For this case, the following chain of equivalences hold
\begin{flushleft}
$T^{(2)}(M)$ is a second-order path term
\allowdisplaybreaks
\begin{align*}
\Leftrightarrow\quad & M\circ^{0\mathbf{T}_{\Sigma^{\boldsymbol{\mathcal{A}}^{(2)}}}(X)}_{s}\mathrm{CH}^{(2)}_{s}\left(\mathfrak{P}^{(2)}\right)\mbox{ is a second-order path term}
\tag{1}
\\
\Leftrightarrow\quad &
\mathrm{ip}^{(2,X)@}_{s}\left(
M\circ^{0\mathbf{T}_{\Sigma^{\boldsymbol{\mathcal{A}}^{(2)}}}(X)}_{s}\mathrm{CH}^{(2)}_{s}\left(\mathfrak{P}^{(2)}\right)
\right)
\mbox{ is a second-order path}
\tag{2}
\\
\Leftrightarrow\quad &
\mathrm{ip}^{(2,X)@}_{s}\left(
M\right)\circ^{0\mathbf{Pth}_{\boldsymbol{\mathcal{A}}^{(2)}}}_{s}
\mathrm{ip}^{(2,X)@}_{s}\left(
\mathrm{CH}^{(2)}_{s}\left(\mathfrak{P}^{(2)}\right)
\right)
\mbox{ is a second-order path}
\tag{3}
\\
\Leftrightarrow\quad &
\mathrm{sc}^{(0,2)}_{s}\left(
\mathrm{ip}^{(2,X)@}_{s}\left(
M\right)\right)
=
\mathrm{tg}^{(0,2)}_{s}\left(
\mathrm{ip}^{(2,X)@}_{s}\left(
\mathrm{CH}^{(2)}_{s}\left(\mathfrak{P}^{(2)}\right)
\right)\right)
\tag{4}
\\
\Leftrightarrow\quad &
\mathrm{sc}^{(0,2)}_{s}\left(
\mathrm{ip}^{(2,X)@}_{s}\left(
N\right)\right)
=
\mathrm{tg}^{(0,2)}_{s}\left(
\mathrm{ip}^{(2,X)@}_{s}\left(
\mathrm{CH}^{(2)}_{s}\left(\mathfrak{P}^{(2)}\right)
\right)\right)
\tag{5}
\\
\Leftrightarrow\quad &
\mathrm{ip}^{(2,X)@}_{s}\left(
N\right)\circ^{0\mathbf{Pth}_{\boldsymbol{\mathcal{A}}^{(2)}}}_{s}
\mathrm{ip}^{(2,X)@}_{s}\left(
\mathrm{CH}^{(2)}_{s}\left(\mathfrak{P}^{(2)}\right)
\right)
\mbox{ is a second-order path}
\tag{6}
\\
\Leftrightarrow\quad &
\mathrm{ip}^{(2,X)@}_{s}\left(
N\circ^{0\mathbf{T}_{\Sigma^{\boldsymbol{\mathcal{A}}^{(2)}}}(X)}_{s}\mathrm{CH}^{(2)}_{s}\left(\mathfrak{P}^{(2)}\right)
\right)
\mbox{ is a second-order path}
\tag{7}
\\
\Leftrightarrow\quad & N\circ^{0\mathbf{T}_{\Sigma^{\boldsymbol{\mathcal{A}}^{(2)}}}(X)}_{s}\mathrm{CH}^{(2)}_{s}\left(\mathfrak{P}^{(2)}\right)\mbox{ is a second-order path term}
\tag{8}
\\
\Leftrightarrow\quad &
T^{(2)}(N) \mbox{ is a second-order path term}.
\tag{9}
\end{align*}
\end{flushleft}

In the just stated chain of equivalences, the first equivalence follows from the description of the second-order elementary translation $T^{(2)}$; the second equivalence follows from Proposition~\ref{PDPT}; 	the third equivalence follows from the fact that $\mathrm{ip}^{(2,X)@}$ is a $\Sigma^{\boldsymbol{\mathcal{A}}^{(2)}}$-homomorphism, according to Definition~\ref{DDIp}; the fourth equivalence follows from the definition of the $0$-composition operation in $\mathbf{Pth}_{\boldsymbol{\mathcal{A}}^{(2)}}$ according to Proposition~\ref{PDPthCatAlg}; the  fifth equivalence follows from the fact that, by hypothesis, $
\mathrm{sc}^{([1],2)}_{t}(
\mathrm{ip}^{(2,X)@}_{t}(
M))
=
\mathrm{sc}^{([1],2)}_{t}(
\mathrm{ip}^{(2,X)@}_{t}(
N))$. In particular, this implies, according to Definition~\ref{DDScTgZ}, that $
\mathrm{sc}^{(0,2)}_{t}(
\mathrm{ip}^{(2,X)@}_{t}(
M))
=
\mathrm{sc}^{(0,2)}_{t}(
\mathrm{ip}^{(2,X)@}_{t}(
N))$; the sixth equivalence follows from the definition of the $0$-composition operation in $\mathbf{Pth}_{\boldsymbol{\mathcal{A}}^{(2)}}$ according to Proposition~\ref{PDPthCatAlg}; the seventh equivalence follows from the fact that $\mathrm{ip}^{(2,X)@}$ is a $\Sigma^{\boldsymbol{\mathcal{A}}}$-homomorphism, according to Definition~\ref{DDIp}; the eighth equivalence follows from Proposition~\ref{PDPT}; finally, the last equivalence follows from the description of the second-order elementary translation $T^{(2)}$.

Assume that either $T^{(2)}(M)$ or $T^{(2)}(N)$ is a second-order path term.
Regarding the $([1],2)$-source, the following chain of equalities holds
\begin{flushleft}
$\mathrm{sc}^{([1],2)}_{s}\left(\mathrm{ip}^{(2,X)@}_{s}\left(
T^{(2)}(M)\right)\right)$
\allowdisplaybreaks
\begin{align*}
&=\mathrm{sc}^{([1],2)}_{s}\left(\mathrm{ip}^{(2,X)@}_{s}\left(
M\circ^{0\mathbf{T}_{\Sigma^{\boldsymbol{\mathcal{A}}^{(2)}}}(X)}_{s}\mathrm{CH}^{(2)}_{s}\left(\mathfrak{P}^{(2)}\right)\right)\right)
\tag{1}
\\&=
\mathrm{sc}^{([1],2)}_{s}\left(
\mathrm{ip}^{(2,X)@}_{s}\left(
M\right)
\circ^{0\mathbf{Pth}_{\boldsymbol{\mathcal{A}}^{(2)}}}_{s}
\mathrm{ip}^{(2,X)@}_{s}\left(
\mathrm{CH}^{(2)}_{s}\left(
\mathfrak{P}^{(2)}\right)
\right)\right)
\tag{2}
\\&=
\mathrm{sc}^{([1],2)}_{s}\left(
\mathrm{ip}^{(2,X)@}_{s}\left(
M\right)
\right)
\circ^{0[\mathbf{PT}_{\boldsymbol{\mathcal{A}}}]}_{s}
\mathrm{sc}^{([1],2)}_{s}\left(
\mathrm{ip}^{(2,X)@}_{s}\left(
\mathrm{CH}^{(2)}_{s}\left(
\mathfrak{P}^{(2)}\right)
\right)\right)
\tag{3}
\\&=
\mathrm{sc}^{([1],2)}_{s}\left(
\mathrm{ip}^{(2,X)@}_{s}\left(
N\right)
\right)
\circ^{0[\mathbf{PT}_{\boldsymbol{\mathcal{A}}}]}_{s}
\mathrm{sc}^{([1],2)}_{s}\left(
\mathrm{ip}^{(2,X)@}_{s}\left(
\mathrm{CH}^{(2)}_{s}\left(
\mathfrak{P}^{(2)}\right)
\right)\right)
\tag{4}
\\&=
\mathrm{sc}^{([1],2)}_{s}\left(
\mathrm{ip}^{(2,X)@}_{s}\left(
N\right)
\circ^{0\mathbf{Pth}_{\boldsymbol{\mathcal{A}}^{(2)}}}_{s}
\mathrm{ip}^{(2,X)@}_{s}\left(
\mathrm{CH}^{(2)}_{s}\left(
\mathfrak{P}^{(2)}\right)
\right)\right)
\tag{5}
\\&=\mathrm{sc}^{([1],2)}_{s}\left(\mathrm{ip}^{(2,X)@}_{s}\left(
N\circ^{0\mathbf{T}_{\Sigma^{\boldsymbol{\mathcal{A}}^{(2)}}}(X)}_{s}\mathrm{CH}^{(2)}_{s}\left(\mathfrak{P}^{(2)}\right)\right)\right)
\tag{6}
\\&=
\mathrm{sc}^{([1],2)}_{s}\left(\mathrm{ip}^{(2,X)@}_{s}\left(
T^{(2)}(N)\right)\right).
\tag{7}
\end{align*}
\end{flushleft}

In the just stated chain of equalities, the first equality follows from the description of the second-order elementary translation $T^{(2)}$; the second equality follows from the fact that $\mathrm{ip}^{(2,X)@}$ is a $\Sigma^{\boldsymbol{\mathcal{A}}^{(2)}}$-homomorphism, according to Definition~\ref{DDIp}; the third equality follows from the fact that  $\mathrm{sc}^{([1],2)}$ is a $\Sigma^{\boldsymbol{\mathcal{A}}}$-homomorphism, according to Proposition~\ref{PDUCatHom}; the fourth equality follows from the fact that, by hypothesis, $
\mathrm{sc}^{([1],2)}_{t}(
\mathrm{ip}^{(2,X)@}_{t}(
M))
=
\mathrm{sc}^{([1],2)}_{t}(
\mathrm{ip}^{(2,X)@}_{t}(
N))$; the fifth equality follows from the fact that  $\mathrm{sc}^{([1],2)}$ is a $\Sigma^{\boldsymbol{\mathcal{A}}}$-homomorphism, according to Proposition~\ref{PDUCatHom};  the sixth equality follows from the fact that $\mathrm{ip}^{(2,X)@}$ is a $\Sigma^{\boldsymbol{\mathcal{A}}^{(2)}}$-homomorphism, according to Definition~\ref{DDIp}; finally, the last equivalence follows from the description of the second-order elementary translation $T^{(2)}$.

The case of the $([1],2)$-target follows by a similar argument.

The remaining case, i.e., the case for which $s=t$ and there is a second-order path $\mathfrak{P}^{(2)}\in\mathrm{Pth}_{\boldsymbol{\mathcal{A}}^{(2)},s}$ such that
\[
T^{(2)}=\mathrm{CH}^{(2)}_{s}\left(\mathfrak{P}^{(2)}\right)\circ^{0\mathbf{T}_{\Sigma^{\boldsymbol{\mathcal{A}}^{(2)}}}(X)}_{s}\underline{\quad},
\]
follows by a similar argument to the proof presented above.

Now, consider the case for the second-order elementary translation $T^{(2)}$ for which  $s=t$ and there is a second-order path $\mathfrak{P}^{(2)}\in\mathrm{Pth}_{\boldsymbol{\mathcal{A}}^{(2)},s}$ such that
\[
T^{(2)}=\underline{\quad}\circ^{1\mathbf{T}_{\Sigma^{\boldsymbol{\mathcal{A}}^{(2)}}}(X)}_{s}\mathrm{CH}^{(2)}_{s}\left(\mathfrak{P}^{(2)}\right).
\]

For this case, the following chain of equivalences hold
\begin{flushleft}
$T^{(2)}(M)$ is a second-order path term
\allowdisplaybreaks
\begin{align*}
\Leftrightarrow\quad & M\circ^{1\mathbf{T}_{\Sigma^{\boldsymbol{\mathcal{A}}^{(2)}}}(X)}_{s}\mathrm{CH}^{(2)}_{s}\left(\mathfrak{P}^{(2)}\right)\mbox{ is a second-order path term}
\tag{1}
\\
\Leftrightarrow\quad &
\mathrm{ip}^{(2,X)@}_{s}\left(
M\circ^{1\mathbf{T}_{\Sigma^{\boldsymbol{\mathcal{A}}^{(2)}}}(X)}_{s}\mathrm{CH}^{(2)}_{s}\left(\mathfrak{P}^{(2)}\right)
\right)
\mbox{ is a second-order path}
\tag{2}
\\
\Leftrightarrow\quad &
\mathrm{ip}^{(2,X)@}_{s}\left(
M\right)\circ^{1\mathbf{Pth}_{\boldsymbol{\mathcal{A}}^{(2)}}}_{s}
\mathrm{ip}^{(2,X)@}_{s}\left(
\mathrm{CH}^{(2)}_{s}\left(\mathfrak{P}^{(2)}\right)
\right)
\mbox{ is a second-order path}
\tag{3}
\\
\Leftrightarrow\quad &
\mathrm{sc}^{([1],2)}_{s}\left(
\mathrm{ip}^{(2,X)@}_{s}\left(
M\right)\right)
=
\mathrm{tg}^{([1],2)}_{s}\left(
\mathrm{ip}^{(2,X)@}_{s}\left(
\mathrm{CH}^{(2)}_{s}\left(\mathfrak{P}^{(2)}\right)
\right)\right)
\tag{4}
\\
\Leftrightarrow\quad &
\mathrm{sc}^{([1],2)}_{s}\left(
\mathrm{ip}^{(2,X)@}_{s}\left(
N\right)\right)
=
\mathrm{tg}^{([1],2)}_{s}\left(
\mathrm{ip}^{(2,X)@}_{s}\left(
\mathrm{CH}^{(2)}_{s}\left(\mathfrak{P}^{(2)}\right)
\right)\right)
\tag{5}
\\
\Leftrightarrow\quad &
\mathrm{ip}^{(2,X)@}_{s}\left(
N\right)\circ^{1\mathbf{Pth}_{\boldsymbol{\mathcal{A}}^{(2)}}}_{s}
\mathrm{ip}^{(2,X)@}_{s}\left(
\mathrm{CH}^{(2)}_{s}\left(\mathfrak{P}^{(2)}\right)
\right)
\mbox{ is a second-order path}
\tag{6}
\\
\Leftrightarrow\quad &
\mathrm{ip}^{(2,X)@}_{s}\left(
N\circ^{1\mathbf{T}_{\Sigma^{\boldsymbol{\mathcal{A}}^{(2)}}}(X)}_{s}\mathrm{CH}^{(2)}_{s}\left(\mathfrak{P}^{(2)}\right)
\right)
\mbox{ is a second-order path}
\tag{7}
\\
\Leftrightarrow\quad & N\circ^{1\mathbf{T}_{\Sigma^{\boldsymbol{\mathcal{A}}^{(2)}}}(X)}_{s}\mathrm{CH}^{(2)}_{s}\left(\mathfrak{P}^{(2)}\right)\mbox{ is a second-order path term}
\tag{8}
\\
\Leftrightarrow\quad &
T^{(2)}(N) \mbox{ is a second-order path term}.
\tag{9}
\end{align*}
\end{flushleft}

In the just stated chain of equivalences, the first equivalence follows from the description of the second-order elementary translation $T^{(2)}$; the second equivalence follows from Proposition~\ref{PDPT}; 	the third equivalence follows from the fact that $\mathrm{ip}^{(2,X)@}$ is a $\Sigma^{\boldsymbol{\mathcal{A}}^{(2)}}$-homomorphism, according to Definition~\ref{DDIp}; the fourth equivalence follows from the definition of the $1$-composition operation in $\mathbf{Pth}_{\boldsymbol{\mathcal{A}}^{(2)}}$ according to Proposition~\ref{PDPthDCatAlg}; the  fifth equivalence follows from the fact that, by hypothesis, $
\mathrm{sc}^{([1],2)}_{t}(
\mathrm{ip}^{(2,X)@}_{t}(
M))
=
\mathrm{sc}^{([1],2)}_{t}(
\mathrm{ip}^{(2,X)@}_{t}(
N))$; the sixth equivalence follows from the definition of the $1$-composition operation in $\mathbf{Pth}_{\boldsymbol{\mathcal{A}}^{(2)}}$ according to Proposition~\ref{PDPthDCatAlg}; the seventh equivalence follows from the fact that $\mathrm{ip}^{(2,X)@}$ is a $\Sigma^{\boldsymbol{\mathcal{A}}^{(2)}}$-homomorphism, according to Definition~\ref{DDIp}; the eighth equivalence follows from Proposition~\ref{PDPT}; finally, the last equivalence follows from the description of the second-order elementary translation $T^{(2)}$.

Assume that either $T^{(2)}(M)$ or $T^{(2)}(N)$ is a secon-order path term.
Regarding the $([1],2)$-source, the following chain of equalities holds
\begin{flushleft}
$\mathrm{sc}^{([1],2)}_{s}\left(\mathrm{ip}^{(2,X)@}_{s}\left(
T^{(2)}(M)\right)\right)$
\allowdisplaybreaks
\begin{align*}
&=\mathrm{sc}^{([1],2)}_{s}\left(\mathrm{ip}^{(2,X)@}_{s}\left(
M\circ^{1\mathbf{T}_{\Sigma^{\boldsymbol{\mathcal{A}}^{(2)}}}(X)}_{s}\mathrm{CH}^{(2)}_{s}\left(\mathfrak{P}^{(2)}\right)\right)\right)
\tag{1}
\\&=
\mathrm{sc}^{([1],2)}_{s}\left(
\mathrm{ip}^{(2,X)@}_{s}\left(
M\right)
\circ^{1\mathbf{Pth}_{\boldsymbol{\mathcal{A}}^{(2)}}}_{s}
\mathrm{ip}^{(2,X)@}_{s}\left(
\mathrm{CH}^{(2)}_{s}\left(
\mathfrak{P}^{(2)}\right)
\right)\right)
\tag{2}
\\&=
\mathrm{sc}^{([1],2)}_{s}\left(
\mathrm{ip}^{(2,X)@}_{s}\left(
\mathrm{CH}^{(2)}_{s}\left(
\mathfrak{P}^{(2)}\right)
\right)\right)
\tag{3}
\\&=
\mathrm{sc}^{([1],2)}_{s}\left(
\mathrm{ip}^{(2,X)@}_{s}\left(
N\right)
\circ^{0\mathbf{Pth}_{\boldsymbol{\mathcal{A}}^{(2)}}}_{s}
\mathrm{ip}^{(2,X)@}_{s}\left(
\mathrm{CH}^{(2)}_{s}\left(
\mathfrak{P}^{(2)}\right)
\right)\right)
\tag{4}
\\&=\mathrm{sc}^{([1],2)}_{s}\left(\mathrm{ip}^{(2,X)@}_{s}\left(
N\circ^{1\mathbf{T}_{\Sigma^{\boldsymbol{\mathcal{A}}^{(2)}}}(X)}_{s}\mathrm{CH}^{(2)}_{s}\left(\mathfrak{P}^{(2)}\right)\right)\right)
\tag{5}
\\&=
\mathrm{sc}^{([1],2)}_{s}\left(\mathrm{ip}^{(2,X)@}_{s}\left(
T^{(2)}(N)\right)\right).
\tag{6}
\end{align*}
\end{flushleft}

In the just stated chain of equalities, the first equality follows from the description of the second-order elementary translation $T^{(2)}$; the second equality follows from the fact that $\mathrm{ip}^{(2,X)@}$ is a $\Sigma^{\boldsymbol{\mathcal{A}}^{(2)}}$-homomorphism, according to Definition~\ref{DDIp}; the third equality follows from the definition of the $1$-composition operation in $\mathbf{Pth}_{\boldsymbol{\mathcal{A}}^{(2)}}$ according to Proposition~\ref{PDPthDCatAlg}; the fourth equality follows from the definition of the $1$-composition operation in $\mathbf{Pth}_{\boldsymbol{\mathcal{A}}^{(2)}}$ according to Proposition~\ref{PDPthDCatAlg}; the fifth equality follows from the fact that $\mathrm{ip}^{(2,X)@}$ is a $\Sigma^{\boldsymbol{\mathcal{A}}^{(2)}}$-homomorphism, according to Definition~\ref{DDIp}; finally, the last equivalence follows from the description of the second-order elementary translation $T^{(2)}$.

Now, regarding the $([1],2)$-target, the following chain of equalities holds
\begin{flushleft}
$\mathrm{tg}^{([1],2)}_{s}\left(\mathrm{ip}^{(2,X)@}_{s}\left(
T^{(2)}(M)\right)\right)$
\allowdisplaybreaks
\begin{align*}
&=\mathrm{tg}^{([1],2)}_{s}\left(\mathrm{ip}^{(2,X)@}_{s}\left(
M\circ^{1\mathbf{T}_{\Sigma^{\boldsymbol{\mathcal{A}}^{(2)}}}(X)}_{s}\mathrm{CH}^{(2)}_{s}\left(\mathfrak{P}^{(2)}\right)\right)\right)
\tag{1}
\\&=
\mathrm{tg}^{([1],2)}_{s}\left(
\mathrm{ip}^{(2,X)@}_{s}\left(
M\right)
\circ^{1\mathbf{Pth}_{\boldsymbol{\mathcal{A}}^{(2)}}}_{s}
\mathrm{ip}^{(2,X)@}_{s}\left(
\mathrm{CH}^{(2)}_{s}\left(
\mathfrak{P}^{(2)}\right)
\right)\right)
\tag{2}
\\&=
\mathrm{tg}^{([1],2)}_{s}\left(
\mathrm{ip}^{(2,X)@}_{s}\left(
M\right)
\right)
\tag{3}
\\&=
\mathrm{tg}^{([1],2)}_{s}\left(
\mathrm{ip}^{(2,X)@}_{s}\left(
N\right)
\right)
\tag{4}
\\&=
\mathrm{tg}^{([1],2)}_{s}\left(
\mathrm{ip}^{(2,X)@}_{s}\left(
N\right)
\circ^{1\mathbf{Pth}_{\boldsymbol{\mathcal{A}}^{(2)}}}_{s}
\mathrm{ip}^{(2,X)@}_{s}\left(
\mathrm{CH}^{(2)}_{s}\left(
\mathfrak{P}^{(2)}\right)
\right)\right)
\tag{5}
\\&=\mathrm{tg}^{([1],2)}_{s}\left(\mathrm{ip}^{(2,X)@}_{s}\left(
N\circ^{1\mathbf{T}_{\Sigma^{\boldsymbol{\mathcal{A}}^{(2)}}}(X)}_{s}\mathrm{CH}^{(2)}_{s}\left(\mathfrak{P}^{(2)}\right)\right)\right)
\tag{6}
\\&=
\mathrm{tg}^{([1],2)}_{s}\left(\mathrm{ip}^{(2,X)@}_{s}\left(
T^{(2)}(N)\right)\right).
\tag{7}
\end{align*}
\end{flushleft}

In the just stated chain of equalities, the first equality follows from the description of the second-order elementary translation $T^{(2)}$; the second equality follows from the fact that $\mathrm{ip}^{(2,X)@}$ is a $\Sigma^{\boldsymbol{\mathcal{A}}^{(2)}}$-homomorphism, according to Definition~\ref{DDIp}; the third equality follows from the definition of the $1$-composition operation in $\mathbf{Pth}_{\boldsymbol{\mathcal{A}}^{(2)}}$ according to Proposition~\ref{PDPthDCatAlg}; the fourth
equality follows from the fact that, by hypothesis, $
\mathrm{tg}^{([1],2)}_{t}(
\mathrm{ip}^{(2,X)@}_{t}(
M))
=
\mathrm{tg}^{([1],2)}_{t}(
\mathrm{ip}^{(2,X)@}_{t}(
N))$; the fifth equality follows from the definition of the $1$-composition operation in $\mathbf{Pth}_{\boldsymbol{\mathcal{A}}^{(2)}}$ according to Proposition~\ref{PDPthDCatAlg}; the sixth equality follows from the fact that $\mathrm{ip}^{(2,X)@}$ is a $\Sigma^{\boldsymbol{\mathcal{A}}^{(2)}}$-homomorphism, according to Definition~\ref{DDIp}; finally, the last equivalence follows from the description of the second-order elementary translation $T^{(2)}$.

The remaining case, i.e., the case for which $s=t$ and there is a second-order path $\mathfrak{P}^{(2)}\in\mathrm{Pth}_{\boldsymbol{\mathcal{A}}^{(2)},s}$ such that
\[
T^{(2)}=\mathrm{CH}^{(2)}_{s}\left(\mathfrak{P}^{(2)}\right)\circ^{1\mathbf{T}_{\Sigma^{\boldsymbol{\mathcal{A}}^{(2)}}}(X)}_{s}\underline{\quad},
\]
follows by a similar argument to the proof presented above.

This complete the base case.

\textsf{Inductive step of the induction.}

Let us suppose that the statement holds for second-order translations of height $m\in \mathbb{N}-\{0\}$. That is, if $T^{(2)}$ is a second-order translation in 
$\mathrm{Tl}_{t}(\mathrm{T}_{\Sigma^{\boldsymbol{\mathcal{A}}^{(2)}}}(X))_{s}$ of height $m$ and $M$ and $N$ are second-order path terms in
$\mathrm{PT}_{\boldsymbol{\mathcal{A}}^{(2)},t}$ satisfying that 
 $$
\left(
\mathrm{ip}^{(2,X)@}_{t}\left(
M
\right), 
\mathrm{ip}^{(2,X)@}_{t}\left(
N
\right)\right)
\in\mathrm{Ker}\left(
\mathrm{sc}^{([1],2)}
\right)_{t}
\cap\mathrm{Ker}\left(
\mathrm{tg}^{([1],2)}
\right)_{t},
$$
then the following properties hold 
\begin{itemize}
\item[(i)]
$T^{(2)}(M)$ is a second-order path term in $\mathrm{PT}_{\boldsymbol{\mathcal{A}}^{(2)},s}$ if, and only if, $T^{(2)}(N)$ is a second-order path term in  $\mathrm{PT}_{\boldsymbol{\mathcal{A}}^{(2)},s}$;
\item[(ii)] If either $T^{(2)}(M)$ or $T^{(2)}(N)$ is a second-order path term in $\mathrm{PT}_{\boldsymbol{\mathcal{A}}^{(2)},s}$, then 
\[
\left(
\mathrm{ip}^{(2,X)@}_{s}\left(
T^{(2)}(M)\right),
\mathrm{ip}^{(2,X)@}_{s}\left(
T^{(2)}(N)\right)
\right)
\in\mathrm{Ker}\left(
\mathrm{sc}^{([1],2)}
\right)_{s}
\cap\mathrm{Ker}\left(
\mathrm{tg}^{([1],2)}
\right)_{s}.
\]
\end{itemize}

Now let $T^{(2)}$ be a second-order translation in 
$\mathrm{Tl}_{t}(\mathrm{T}_{\Sigma^{\boldsymbol{\mathcal{A}}^{(2)}}}(X))_{s}$ of height $m+1$ and let $M$ and $N$ be second-order path terms in
$\mathrm{PT}_{\boldsymbol{\mathcal{A}}^{(2)},t}$ satisfying that 
 $$
\left(
\mathrm{ip}^{(2,X)@}_{t}\left(
M
\right), 
\mathrm{ip}^{(2,X)@}_{t}\left(
N
\right)\right)
\in\mathrm{Ker}\left(
\mathrm{sc}^{([1],2)}
\right)_{t}
\cap\mathrm{Ker}\left(
\mathrm{tg}^{([1],2)}
\right)_{t}.
$$

We consider the different possibilities for the case of $T^{(2)}$ being a second-order  translation according to Definition~\ref{DDTrans}.

Assume that there is a word $\mathbf{s}\in S^{\star}-\{\lambda\}$, an index $k\in\bb{\mathbf{s}}$, an operation symbol $\sigma\in \Sigma_{\mathbf{s},s}$, a family of second-order paths $(\mathfrak{P}^{(2)}_{j})_{j\in k}\in \prod_{j\in k}\mathrm{Pth}_{\boldsymbol{\mathcal{A}}^{(2)},s_{j}}$, a family of second-order paths $(\mathfrak{P}^{(2)}_{l})_{l\in \bb{\mathbf{s}}-(k+1)}\in \prod_{l\in \bb{\mathbf{s}}-(k+1)}\mathrm{Pth}_{\boldsymbol{\mathcal{A}}^{(2)},s_{l}}$ 
and a second-order translation $T^{(1)'}$ in $\mathrm{Tl}_{t}(\mathrm{T}_{\Sigma^{\boldsymbol{\mathcal{A}}^{(2)}}}(X))_{s_{k}}$ of height $m$  such that 
\begin{multline*}
T^{(2)}=\sigma^{\mathbf{T}_{\Sigma^{\boldsymbol{\mathcal{A}}^{(2)}}}(X)}\left(
\mathrm{CH}^{(2)}_{s_{0}}\left(\mathfrak{P}^{(2)}_{0}\right),
\cdots,
\mathrm{CH}^{(2)}_{s_{k-1}}\left(\mathfrak{P}^{(2)}_{k-1}\right),
\right.
\\
\left.
T^{(2)'},
\mathrm{CH}^{(2)}_{s_{k+1}}\left(\mathfrak{P}^{(2)}_{k+1}\right),
\cdots,
\mathrm{CH}^{(2)}_{s_{\bb{\mathbf{s}}-1}}\left(\mathfrak{P}^{(2)}_{\bb{\mathbf{s}}-1}\right)
\right).
\end{multline*}

Note that, for this case, the following equivalences hold
\begin{flushleft}
$T^{(2)}(M)$ is a second-order path term
\allowdisplaybreaks
\begin{align*}
\Leftrightarrow\quad&
\sigma^{\mathbf{T}_{\Sigma^{\boldsymbol{\mathcal{A}}^{(2)}}}(X)}\left(
\mathrm{CH}^{(2)}_{s_{0}}\left(\mathfrak{P}^{(2)}_{0}\right),
\cdots,
\mathrm{CH}^{(2)}_{s_{k-1}}\left(\mathfrak{P}^{(2)}_{k-1}\right),
T^{(2)'}(M),
\right.
\\&\qquad\qquad\qquad\qquad\quad
\left.
\mathrm{CH}^{(2)}_{s_{k+1}}\left(\mathfrak{P}^{(2)}_{k+1}\right),
\cdots,
\mathrm{CH}^{(2)}_{s_{\bb{\mathbf{s}}-1}}\left(\mathfrak{P}^{(2)}_{\bb{\mathbf{s}}-1}\right)
\right)
\\&\qquad\qquad\qquad\qquad\qquad\qquad\qquad\qquad\qquad
\mbox{ is a second-order path term}
\tag{1}
\\
\Leftrightarrow\quad&
\mathrm{ip}^{(2,X)@}_{s}\left(
\sigma^{\mathbf{T}_{\Sigma^{\boldsymbol{\mathcal{A}}^{(2)}}}(X)}\left(
\mathrm{CH}^{(2)}_{s_{0}}\left(\mathfrak{P}^{(2)}_{0}\right),
\cdots,
\mathrm{CH}^{(2)}_{s_{k-1}}\left(\mathfrak{P}^{(2)}_{k-1}\right),
T^{(2)'}(M),
\right.\right.
\\&\qquad\qquad\qquad\qquad\quad
\left.\left.
\mathrm{CH}^{(2)}_{s_{k+1}}\left(\mathfrak{P}^{(2)}_{k+1}\right),
\cdots,
\mathrm{CH}^{(2)}_{s_{\bb{\mathbf{s}}-1}}\left(\mathfrak{P}^{(2)}_{\bb{\mathbf{s}}-1}\right)
\right)\right)
\\&\qquad\qquad\qquad\qquad\qquad\qquad\qquad\qquad\qquad\qquad
\mbox{ is a second-order path}
\tag{2}
\\
\Leftrightarrow\quad&
\sigma^{\mathbf{Pth}_{\boldsymbol{\mathcal{A}}^{(2)}}}\left(
\mathrm{ip}^{(2,X)@}_{s_{0}}\left(
\mathrm{CH}^{(2)}_{s_{0}}\left(\mathfrak{P}^{(2)}_{0}\right)\right),
\cdots,
\mathrm{ip}^{(2,X)@}_{s_{k-1}}\left(
\mathrm{CH}^{(2)}_{s_{k-1}}\left(\mathfrak{P}^{(2)}_{k-1}\right)\right),
\right.
\\&\qquad\qquad\qquad
\mathrm{ip}^{(2,X)@}_{s_{k}}\left(
T^{(2)'}(M)\right),
\mathrm{ip}^{(2,X)@}_{s_{k+1}}\left(
\mathrm{CH}^{(2)}_{s_{k+1}}\left(\mathfrak{P}^{(2)}_{k+1}\right)\right),
\cdots,
\\&\qquad\qquad\qquad\qquad\qquad\qquad\qquad
\left.
\mathrm{ip}^{(2,X)@}_{s_{\bb{\mathbf{s}}-1}}\left(
\mathrm{CH}^{(2)}_{s_{\bb{\mathbf{s}}-1}}\left(\mathfrak{P}^{(2)}_{\bb{\mathbf{s}}-1}\right)
\right)\right)
\\&\qquad\qquad\qquad\qquad\qquad\qquad\qquad\qquad\qquad\qquad
\mbox{ is a second-order path}
\tag{3}
\\
\Leftrightarrow\quad&
\mathrm{ip}^{(2,X)@}_{s_{k}}\left(
T^{(2)'}(M)\right)
\mbox{ is a secod-order path}
\tag{4}
\\
\Leftrightarrow\quad&
T^{(2)'}(M)
\mbox{ is a second-order path term}
\tag{5}
\\
\Leftrightarrow\quad&
T^{(2)'}(N)
\mbox{ is a second-order path term}
\tag{6}
\\
\Leftrightarrow\quad&
\mathrm{ip}^{(2,X)@}_{s_{k}}\left(
T^{(2)'}(N)\right)
\mbox{ is a second-order path}
\tag{7}
\\
\Leftrightarrow\quad&
\sigma^{\mathbf{Pth}_{\boldsymbol{\mathcal{A}}^{(2)}}}\left(
\mathrm{ip}^{(2,X)@}_{s_{0}}\left(
\mathrm{CH}^{(2)}_{s_{0}}\left(\mathfrak{P}^{(2)}_{0}\right)\right),
\cdots,
\mathrm{ip}^{(2,X)@}_{s_{k-1}}\left(
\mathrm{CH}^{(2)}_{s_{k-1}}\left(\mathfrak{P}^{(2)}_{k-1}\right)\right),
\right.
\\&\qquad\qquad\qquad
\mathrm{ip}^{(2,X)@}_{s_{k}}\left(
T^{(2)'}(N)\right),
\mathrm{ip}^{(2,X)@}_{s_{k+1}}\left(
\mathrm{CH}^{(2)}_{s_{k+1}}\left(\mathfrak{P}^{(2)}_{k+1}\right)\right),
\cdots,
\\&\qquad\qquad\qquad\qquad\qquad\qquad\qquad
\left.
\mathrm{ip}^{(2,X)@}_{s_{\bb{\mathbf{s}}-1}}\left(
\mathrm{CH}^{(2)}_{s_{\bb{\mathbf{s}}-1}}\left(\mathfrak{P}^{(2)}_{\bb{\mathbf{s}}-1}\right)
\right)\right)
\\&\qquad\qquad\qquad\qquad\qquad\qquad\qquad\qquad\qquad\qquad
\mbox{ is a second-order path}
\tag{8}
\\
\Leftrightarrow\quad&
\mathrm{ip}^{(2,X)@}_{s}\left(
\sigma^{\mathbf{T}_{\Sigma^{\boldsymbol{\mathcal{A}}^{(2)}}}(X)}\left(
\mathrm{CH}^{(2)}_{s_{0}}\left(\mathfrak{P}^{(2)}_{0}\right),
\cdots,
\mathrm{CH}^{(2)}_{s_{k-1}}\left(\mathfrak{P}^{(2)}_{k-1}\right),
T^{(2)'}(N),
\right.\right.
\\&\qquad\qquad\qquad\qquad\quad
\left.\left.
\mathrm{CH}^{(2)}_{s_{k+1}}\left(\mathfrak{P}^{(2)}_{k+1}\right),
\cdots,
\mathrm{CH}^{(2)}_{s_{\bb{\mathbf{s}}-1}}\left(\mathfrak{P}^{(2)}_{\bb{\mathbf{s}}-1}\right)
\right)\right)
\\&\qquad\qquad\qquad\qquad\qquad\qquad\qquad\qquad\qquad\qquad
\mbox{ is a second-order path}
\tag{9}
\\
\Leftrightarrow\quad&
\sigma^{\mathbf{T}_{\Sigma^{\boldsymbol{\mathcal{A}}^{(2)}}}(X)}\left(
\mathrm{CH}^{(2)}_{s_{0}}\left(\mathfrak{P}^{(2)}_{0}\right),
\cdots,
\mathrm{CH}^{(2)}_{s_{k-1}}\left(\mathfrak{P}^{(2)}_{k-1}\right),
T^{(2)'}(N),
\right.
\\&\qquad\qquad\qquad\qquad\quad
\left.
\mathrm{CH}^{(2)}_{s_{k+1}}\left(\mathfrak{P}^{(2)}_{k+1}\right),
\cdots,
\mathrm{CH}^{(2)}_{s_{\bb{\mathbf{s}}-1}}\left(\mathfrak{P}^{(2)}_{\bb{\mathbf{s}}-1}\right)
\right)
\\&\qquad\qquad\qquad\qquad\qquad\qquad\qquad\qquad\qquad
\mbox{ is a second-order path term}
\tag{10}
\\
\Leftrightarrow\quad&
T^{(2)}(N) \mbox{ is a second-order path term}.
\tag{11}
\end{align*}
\end{flushleft}

In the just stated chain of equivalences, the first equivalence follows from the description of the second-order translation $T^{(2)}$; the second equivalence follows from Proposition~\ref{PDPT}; 	the third equivalence follows from the fact that $\mathrm{ip}^{(2,X)@}$ is a $\Sigma^{\boldsymbol{\mathcal{A}}^{(2)}}$-homomorphism, according to Definition~\ref{DDIp}; the fourth equivalence follows from the description of the operation $\sigma$ in the many-sorted partial $\Sigma^{\boldsymbol{\mathcal{A}}^{(2)}}$-algebra $\mathbf{Pth}_{\boldsymbol{\mathcal{A}}^{(2)}}$, according to Proposition~\ref{PDPthAlg}; the fifth equivalence follows from Proposition~\ref{PDPT}; the sixth equivalence follows by induction; the seventh equivalence follows from Proposition~\ref{PDPT}; the eighth equivalence follows from the description of the operation $\sigma$ in the many-sorted partial $\Sigma^{\boldsymbol{\mathcal{A}}^{(2)}}$-algebra $\mathbf{Pth}_{\boldsymbol{\mathcal{A}}^{(2)}}$, according to Proposition~\ref{PDPthAlg}; the ninth equivalence follows from the fact that $\mathrm{ip}^{(2,X)@}$ is a $\Sigma^{\boldsymbol{\mathcal{A}}^{(2)}}$-homomorphism, according to Definition~\ref{DDIp}; the tenth equivalence follows from Proposition~\ref{PDPT}; finally, the last equivalence follows from the description of the second-order translation $T^{(2)}$.

Therefore, the first condition of the proposition holds.

Now, regarding the $([1],2)$-source, the following chain of equalities holds
\begin{flushleft}
$\mathrm{sc}^{([1],2)}_{s}\left(\mathrm{ip}^{(2,X)@}_{s}\left(
T^{(2)}(M)\right)\right)$
\allowdisplaybreaks
\begin{align*}
&=
\mathrm{sc}^{([1],2)}_{s}\left(
\mathrm{ip}^{(2,X)@}_{s}\left(
\sigma^{\mathbf{T}_{\Sigma^{\boldsymbol{\mathcal{A}}^{(2)}}}(X)}\left(
\mathrm{CH}^{(2)}_{s_{0}}\left(\mathfrak{P}^{(2)}_{0}\right),
\cdots,
\mathrm{CH}^{(2)}_{s_{k-1}}\left(\mathfrak{P}^{(2)}_{k-1}\right),
\right.\right.\right.
\\&\qquad\qquad\qquad\qquad
\left.\left.\left.
T^{(2)'}(M),
\mathrm{CH}^{(2)}_{s_{k+1}}\left(\mathfrak{P}^{(2)}_{k+1}\right),
\cdots,
\mathrm{CH}^{(2)}_{s_{\bb{\mathbf{s}}-1}}\left(\mathfrak{P}^{(2)}_{\bb{\mathbf{s}}-1}\right)
\right)
\right)\right)
\tag{1}
\\&=
\mathrm{sc}^{([1],2)}_{s}\left(
\sigma^{\mathbf{Pth}_{\boldsymbol{\mathcal{A}}^{(2)}}}\left(
\mathrm{ip}^{(2,X)@}_{s_{0}}\left(
\mathrm{CH}^{(2)}_{s_{0}}\left(\mathfrak{P}^{(2)}_{0}\right)\right),
\cdots,
\mathrm{ip}^{(2,X)@}_{s_{k-1}}\left(
\mathrm{CH}^{(2)}_{s_{k-1}}\left(\mathfrak{P}^{(2)}_{k-1}\right)\right),
\right.\right.
\\&\qquad\qquad\qquad\qquad
\mathrm{ip}^{(2,X)@}_{t}\left(
T^{(2)'}(M)\right),
\mathrm{ip}^{(2,X)@}_{s_{k+1}}\left(
\mathrm{CH}^{(2)}_{s_{k+1}}\left(\mathfrak{P}^{(2)}_{k+1}\right)\right),
\\&\qquad\qquad\qquad\qquad\qquad\qquad\qquad\qquad\quad
\left.\left.
\cdots,
\mathrm{ip}^{(2,X)@}_{s_{\bb{\mathbf{s}}-1}}\left(
\mathrm{CH}^{(2)}_{s_{\bb{\mathbf{s}}-1}}\left(\mathfrak{P}^{(2)}_{\bb{\mathbf{s}}-1}\right)\right)
\right)
\right)
\tag{2}
\\&=
\sigma^{[\mathbf{PT}_{\boldsymbol{\mathcal{A}}}]}\left(
\mathrm{sc}^{([1],2)}_{s_{0}}\left(
\mathrm{ip}^{(2,X)@}_{s_{0}}\left(
\mathrm{CH}^{(2)}_{s_{0}}\left(\mathfrak{P}^{(2)}_{0}\right)\right)\right),
\cdots,
\right.
\\&\qquad\qquad
\mathrm{sc}^{([1],2)}_{s_{k-1}}\left(
\mathrm{ip}^{(2,X)@}_{s_{k-1}}\left(
\mathrm{CH}^{(2)}_{s_{k-1}}\left(\mathfrak{P}^{(2)}_{k-1}\right)\right)\right),
\mathrm{sc}^{([1],2)}_{t}\left(
\mathrm{ip}^{(2,X)@}_{t}\left(
T^{(2)'}(M)\right)\right),
\\&\qquad\qquad\qquad\qquad
\mathrm{sc}^{([1],2)}_{s_{k+1}}\left(
\mathrm{ip}^{(2,X)@}_{s_{k+1}}\left(
\mathrm{CH}^{(2)}_{s_{k+1}}\left(\mathfrak{P}^{(2)}_{k+1}\right)\right)\right),
\\&\qquad\qquad\qquad\qquad\qquad\qquad
\left.
\cdots,
\mathrm{sc}^{([1],2)}_{s_{\bb{\mathbf{s}}-1}}\left(
\mathrm{ip}^{(2,X)@}_{s_{\bb{\mathbf{s}}-1}}\left(
\mathrm{CH}^{(2)}_{s_{\bb{\mathbf{s}}-1}}\left(\mathfrak{P}^{(2)}_{\bb{\mathbf{s}}-1}\right)\right)
\right)
\right)
\tag{3}
\\&=
\sigma^{[\mathbf{PT}_{\boldsymbol{\mathcal{A}}}]}\left(
\mathrm{sc}^{([1],2)}_{s_{0}}\left(
\mathrm{ip}^{(2,X)@}_{s_{0}}\left(
\mathrm{CH}^{(2)}_{s_{0}}\left(\mathfrak{P}^{(2)}_{0}\right)\right)\right),
\cdots,
\right.
\\&\qquad\qquad
\mathrm{sc}^{([1],2)}_{s_{k-1}}\left(
\mathrm{ip}^{(2,X)@}_{s_{k-1}}\left(
\mathrm{CH}^{(2)}_{s_{k-1}}\left(\mathfrak{P}^{(2)}_{k-1}\right)\right)\right),
\mathrm{sc}^{([1],2)}_{t}\left(
\mathrm{ip}^{(2,X)@}_{t}\left(
T^{(2)'}(N)\right)\right),
\\&\qquad\qquad\qquad\qquad
\mathrm{sc}^{([1],2)}_{s_{k+1}}\left(
\mathrm{ip}^{(2,X)@}_{s_{k+1}}\left(
\mathrm{CH}^{(2)}_{s_{k+1}}\left(\mathfrak{P}^{(2)}_{k+1}\right)\right)\right),
\\&\qquad\qquad\qquad\qquad\qquad\qquad
\left.
\cdots,
\mathrm{sc}^{([1],2)}_{s_{\bb{\mathbf{s}}-1}}\left(
\mathrm{ip}^{(2,X)@}_{s_{\bb{\mathbf{s}}-1}}\left(
\mathrm{CH}^{(2)}_{s_{\bb{\mathbf{s}}-1}}\left(\mathfrak{P}^{(2)}_{\bb{\mathbf{s}}-1}\right)\right)
\right)
\right)
\tag{4}
\\&=
\mathrm{sc}^{([1],2)}_{s}\left(
\sigma^{\mathbf{Pth}_{\boldsymbol{\mathcal{A}}^{(2)}}}\left(
\mathrm{ip}^{(2,X)@}_{s_{0}}\left(
\mathrm{CH}^{(2)}_{s_{0}}\left(\mathfrak{P}^{(2)}_{0}\right)\right),
\cdots,
\mathrm{ip}^{(2,X)@}_{s_{k-1}}\left(
\mathrm{CH}^{(2)}_{s_{k-1}}\left(\mathfrak{P}^{(2)}_{k-1}\right)\right),
\right.\right.
\\&\qquad\qquad\qquad\qquad
\mathrm{ip}^{(2,X)@}_{t}\left(
T^{(2)'}(N)\right),
\mathrm{ip}^{(2,X)@}_{s_{k+1}}\left(
\mathrm{CH}^{(2)}_{s_{k+1}}\left(\mathfrak{P}^{(2)}_{k+1}\right)\right),
\\&\qquad\qquad\qquad\qquad\qquad\qquad\qquad\qquad\quad
\left.\left.
\cdots,
\mathrm{ip}^{(2,X)@}_{s_{\bb{\mathbf{s}}-1}}\left(
\mathrm{CH}^{(2)}_{s_{\bb{\mathbf{s}}-1}}\left(\mathfrak{P}^{(2)}_{\bb{\mathbf{s}}-1}\right)\right)
\right)
\right)
\tag{5}
\\&
\mathrm{sc}^{([1],2)}_{s}\left(
\mathrm{ip}^{(2,X)@}_{s}\left(
\sigma^{\mathbf{T}_{\Sigma^{\boldsymbol{\mathcal{A}}^{(2)}}}(X)}\left(
\mathrm{CH}^{(2)}_{s_{0}}\left(\mathfrak{P}^{(2)}_{0}\right),
\cdots,
\mathrm{CH}^{(2)}_{s_{k-1}}\left(\mathfrak{P}^{(2)}_{k-1}\right),
\right.\right.\right.
\\&\qquad\qquad\qquad\qquad
\left.\left.\left.
T^{(2)'}(N),
\mathrm{CH}^{(2)}_{s_{k+1}}\left(\mathfrak{P}^{(2)}_{k+1}\right),
\cdots,
\mathrm{CH}^{(2)}_{s_{\bb{\mathbf{s}}-1}}\left(\mathfrak{P}^{(2)}_{\bb{\mathbf{s}}-1}\right)
\right)
\right)\right)
\tag{6}
\\&=\mathrm{sc}^{([1],2)}_{s}\left(\mathrm{ip}^{(2,X)@}_{s}\left(
T^{(2)}(N)\right)\right).
\tag{7}
\end{align*}
\end{flushleft}

In the just stated chain of equalities, the first equality unravels the description of the  second-order elementary translation $T^{(2)}$; the second equality follows from the fact that $\mathrm{ip}^{(2,X)@}$ is a $\Sigma^{\boldsymbol{\mathcal{A}}^{(2)}}$-homomorphism, according to Definition~\ref{DDIp}; the third equality follows from the fact that $\mathrm{sc}^{([1],2)}$ is a $\Sigma^{\boldsymbol{\mathcal{A}}}$-homomorphism, according to Proposition~\ref{PDUCatHom}; the fourth equality follows from the fact that, by induction, $
\mathrm{sc}^{([1],2)}_{t}(
\mathrm{ip}^{(2,X)@}_{t}(
T^{(2)'}(M)))
=
\mathrm{sc}^{([1],2)}_{t}(
\mathrm{ip}^{(2,X)@}_{t}(
T^{(2)'}(N)))$;
the fifth equality follows from the fact that $\mathrm{sc}^{([1],2)}$ is a $\Sigma^{\boldsymbol{\mathcal{A}}}$-homomorphism, according to Proposition~\ref{PDUCatHom}; the sixth equality follows from the fact that $\mathrm{ip}^{(2,X)@}$ is a $\Sigma^{\boldsymbol{\mathcal{A}}^{(2)}}$-homomorphism, according to Definition~\ref{DDIp}; finally, the last equality recovers  the description of the second-order translation $T^{(2)}$.

The case of the $([1],2)$-target follows by a similar argument.

Now, regarding the other options for the second-order  translation $T^{(2)}$, it could be the case that $s=t$, there is a second-order path $\mathfrak{P}^{(2)}\in\mathrm{Pth}_{\boldsymbol{\mathcal{A}}^{(2)},s}$ and a second-order translation $T^{(2)'}$ in $\mathrm{Tl}_{t}(\mathrm{T}_{\Sigma^{\boldsymbol{\mathcal{A}}^{(2)}}}(X))_{s}$ of height $m$  such that  such that
\[
T^{(2)}=T^{(2)'}\circ^{0\mathbf{T}_{\Sigma^{\boldsymbol{\mathcal{A}}^{(2)}}}(X)}_{s}\mathrm{CH}^{(2)}_{s}\left(\mathfrak{P}^{(2)}\right).
\]

For this case, the following chain of equivalences hold
\begin{flushleft}
$T^{(2)}(M)$ is a second-order path term
\allowdisplaybreaks
\begin{align*}
\Leftrightarrow\quad & T^{(2)'}(M)\circ^{0\mathbf{T}_{\Sigma^{\boldsymbol{\mathcal{A}}^{(2)}}}(X)}_{s}\mathrm{CH}^{(2)}_{s}\left(\mathfrak{P}^{(2)}\right)\mbox{ is a second-order path term}
\tag{1}
\\
\Leftrightarrow\quad &
\mathrm{ip}^{(2,X)@}_{s}\left(
T^{(2)'}(M)\circ^{0\mathbf{T}_{\Sigma^{\boldsymbol{\mathcal{A}}^{(2)}}}(X)}_{s}\mathrm{CH}^{(2)}_{s}\left(\mathfrak{P}^{(2)}\right)
\right)
\mbox{ is a second-order path}
\tag{2}
\\
\Leftrightarrow\quad &
\mathrm{ip}^{(2,X)@}_{s}\left(
T^{(2)'}(M)\right)\circ^{0\mathbf{Pth}_{\boldsymbol{\mathcal{A}}^{(2)}}}_{s}
\mathrm{ip}^{(2,X)@}_{s}\left(
\mathrm{CH}^{(2)}_{s}\left(\mathfrak{P}^{(2)}\right)
\right)
\\&\qquad\qquad\qquad\qquad\qquad\qquad\qquad\qquad\qquad\qquad\quad
\mbox{ is a second-order path}
\tag{3}
\\
\Leftrightarrow\quad &
\begin{cases}
\mathrm{ip}^{(2,X)@}_{s}\left(
T^{(2)'}(M)\right)\mbox{ is a secod-order path}
\\
\mathrm{sc}^{(0,2)}_{s}\left(
\mathrm{ip}^{(2,X)@}_{s}\left(
T^{(2)'}(M)\right)\right)
=
\mathrm{tg}^{(0,2)}_{s}\left(
\mathrm{ip}^{(2,X)@}_{s}\left(
\mathrm{CH}^{(2)}_{s}\left(\mathfrak{P}^{(2)}\right)
\right)\right)
\end{cases}
\tag{4}
\\
\Leftrightarrow\quad &
\begin{cases}
T^{(2)'}(M)\mbox{ is a second-order path term}
\\
\mathrm{sc}^{(0,2)}_{s}\left(
\mathrm{ip}^{(2,X)@}_{s}\left(
T^{(2)'}(M)\right)\right)
=
\mathrm{tg}^{(0,2)}_{s}\left(
\mathrm{ip}^{(2,X)@}_{s}\left(
\mathrm{CH}^{(2)}_{s}\left(\mathfrak{P}^{(2)}\right)
\right)\right)
\end{cases}
\tag{5}
\\
\Leftrightarrow\quad &
\begin{cases}
T^{(2)'}(N)\mbox{ is a second-order path term}
\\
\mathrm{sc}^{(0,2)}_{s}\left(
\mathrm{ip}^{(2,X)@}_{s}\left(
T^{(2)'}(N)\right)\right)
=
\mathrm{tg}^{(0,2)}_{s}\left(
\mathrm{ip}^{(2,X)@}_{s}\left(
\mathrm{CH}^{(2)}_{s}\left(\mathfrak{P}^{(2)}\right)
\right)\right)
\end{cases}
\tag{6}
\\
\Leftrightarrow\quad &
\begin{cases}
\mathrm{ip}^{(2,X)@}_{s}\left(
T^{(2)'}(N)\right)\mbox{ is a second-order path}
\\
\mathrm{sc}^{(0,2)}_{s}\left(
\mathrm{ip}^{(2,X)@}_{s}\left(
T^{(2)'}(N)\right)\right)
=
\mathrm{tg}^{(0,2)}_{s}\left(
\mathrm{ip}^{(2,X)@}_{s}\left(
\mathrm{CH}^{(2)}_{s}\left(\mathfrak{P}^{(2)}\right)
\right)\right)
\end{cases}
\tag{7}
\\
\Leftrightarrow\quad &
\mathrm{ip}^{(2,X)@}_{s}\left(
T^{(2)'}(N)\right)\circ^{0\mathbf{Pth}_{\boldsymbol{\mathcal{A}}^{(2)}}}_{s}
\mathrm{ip}^{(2,X)@}_{s}\left(
\mathrm{CH}^{(2)}_{s}\left(\mathfrak{P}^{(2)}\right)
\right)
\\&\qquad\qquad\qquad\qquad\qquad\qquad\qquad\qquad\qquad\qquad\quad
\mbox{ is a second-order path}
\tag{8}
\\
\Leftrightarrow\quad &
\mathrm{ip}^{(2,X)@}_{s}\left(
T^{(2)'}(N)\circ^{0\mathbf{T}_{\Sigma^{\boldsymbol{\mathcal{A}}^{(2)}}}(X)}_{s}\mathrm{CH}^{(2)}_{s}\left(\mathfrak{P}^{(2)}\right)
\right)
\mbox{ is a second-order path}
\tag{9}
\\
\Leftrightarrow\quad & T^{(2)'}(N)\circ^{0\mathbf{T}_{\Sigma^{\boldsymbol{\mathcal{A}}^{(2)}}}(X)}_{s}\mathrm{CH}^{(2)}_{s}\left(\mathfrak{P}^{(2)}\right)\mbox{ is a second-order path term}
\tag{10}
\\
\Leftrightarrow\quad &
T^{(2)}(N) \mbox{ is a second-order path term}.
\tag{11}
\end{align*}
\end{flushleft}

In the just stated chain of equivalences, the first equivalence follows from the description of the second-order  translation $T^{(2)}$; the second equivalence follows from Proposition~\ref{PDPT}; 	the third equivalence follows from the fact that $\mathrm{ip}^{(2,X)@}$ is a $\Sigma^{\boldsymbol{\mathcal{A}}^{(2)}}$-homomorphism, according to Definition~\ref{DDIp}; the fourth equivalence follows from the definition of the $0$-composition operation in $\mathbf{Pth}_{\boldsymbol{\mathcal{A}}^{(2)}}$ according to Proposition~\ref{PDPthCatAlg}; 
the fifth equivalence follows from Proposition~\ref{PDPT};
the  sixth equivalence follows by induction. Moreover, since the equality $
\mathrm{sc}^{([1],2)}_{t}(
\mathrm{ip}^{(2,X)@}_{t}(
T^{(2)'}(M)))
=
\mathrm{sc}^{([1],2)}_{t}(
\mathrm{ip}^{(2,X)@}_{t}(
T^{(2)'}(N)))$ also holds by induction, the equality $
\mathrm{sc}^{(0,2)}_{t}(
\mathrm{ip}^{(2,X)@}_{t}(
T^{(2)'}(M)))
=
\mathrm{sc}^{(0,2)}_{t}(
\mathrm{ip}^{(2,X)@}_{t}(
T^{(2)'}(N)))$ holds according to Definition~\ref{DDScTgZ}; the seventh equivalence follows from Proposition~\ref{PDPT}; the eighth equivalence follows from the definition of the $0$-composition operation in $\mathbf{Pth}_{\boldsymbol{\mathcal{A}}^{(2)}}$ according to Proposition~\ref{PDPthCatAlg}; the ninth equivalence follows from the fact that $\mathrm{ip}^{(2,X)@}$ is a $\Sigma^{\boldsymbol{\mathcal{A}}^{(2)}}$-homomorphism, according to Definition~\ref{DDIp}; the tenth equivalence follows from Proposition~\ref{PDPT}; finally, the last equivalence follows from the description of the second-order translation $T^{(2)}$.

Assume that either $T^{(2)}(M)$ or $T^{(2)}(N)$ is a second-order path term.
Regarding the $([1],2)$-source, the following chain of equalities holds
\begin{flushleft}
$\mathrm{sc}^{([1],2)}_{s}\left(\mathrm{ip}^{(2,X)@}_{s}\left(
T^{(2)}(M)\right)\right)$
\allowdisplaybreaks
\begin{align*}
&=\mathrm{sc}^{([1],2)}_{s}\left(\mathrm{ip}^{(2,X)@}_{s}\left(
T^{(2)'}(M)\circ^{0\mathbf{T}_{\Sigma^{\boldsymbol{\mathcal{A}}^{(2)}}}(X)}_{s}\mathrm{CH}^{(2)}_{s}\left(\mathfrak{P}^{(2)}\right)\right)\right)
\tag{1}
\\&=
\mathrm{sc}^{([1],2)}_{s}\left(
\mathrm{ip}^{(2,X)@}_{s}\left(
T^{(2)'}(M)\right)
\circ^{0\mathbf{Pth}_{\boldsymbol{\mathcal{A}}^{(2)}}}_{s}
\mathrm{ip}^{(2,X)@}_{s}\left(
\mathrm{CH}^{(2)}_{s}\left(
\mathfrak{P}^{(2)}\right)
\right)\right)
\tag{2}
\\&=
\mathrm{sc}^{([1],2)}_{s}\left(
\mathrm{ip}^{(2,X)@}_{s}\left(
T^{(2)'}(M)\right)\right)
\circ^{0[\mathbf{PT}_{\boldsymbol{\mathcal{A}}}]}_{s}
\mathrm{sc}^{([1],2)}_{s}\left(
\mathrm{ip}^{(2,X)@}_{s}\left(
\mathrm{CH}^{(2)}_{s}\left(
\mathfrak{P}^{(2)}\right)
\right)\right)
\tag{3}
\\&=
\mathrm{sc}^{([1],2)}_{s}\left(
\mathrm{ip}^{(2,X)@}_{s}\left(
T^{(2)'}(N)\right)\right)
\circ^{0[\mathbf{PT}_{\boldsymbol{\mathcal{A}}}]}_{s}
\mathrm{sc}^{([1],2)}_{s}\left(
\mathrm{ip}^{(2,X)@}_{s}\left(
\mathrm{CH}^{(2)}_{s}\left(
\mathfrak{P}^{(2)}\right)
\right)\right)
\tag{4}
\\&=
\mathrm{sc}^{([1],2)}_{s}\left(
\mathrm{ip}^{(1,X)@}_{s}\left(
T^{(2)'}(N)\right)
\circ^{0\mathbf{Pth}_{\boldsymbol{\mathcal{A}}^{(2)}}}_{s}
\mathrm{ip}^{(2,X)@}_{s}\left(
\mathrm{CH}^{(2)}_{s}\left(
\mathfrak{P}^{(2)}\right)
\right)\right)
\tag{5}
\\&=\mathrm{sc}^{([1],2)}_{s}\left(\mathrm{ip}^{(2,X)@}_{s}\left(
T^{(1)'}(N)\circ^{0\mathbf{T}_{\Sigma^{\boldsymbol{\mathcal{A}}^{(2)}}}(X)}_{s}\mathrm{CH}^{(2)}_{s}\left(\mathfrak{P}^{(2)}\right)\right)\right)
\tag{6}
\\&=
\mathrm{sc}^{([1],2)}_{s}\left(\mathrm{ip}^{(2,X)@}_{s}\left(
T^{(2)}(N)\right)\right).
\tag{7}
\end{align*}
\end{flushleft}

In the just stated chain of equalities, the first equality follows from the description of the second-order  translation $T^{(2)}$; the second equality follows from the fact that $\mathrm{ip}^{(2,X)@}$ is a $\Sigma^{\boldsymbol{\mathcal{A}}^{(2)}}$-homomorphism, according to Definition~\ref{DDIp}; the third equality follows from the fact that  $\mathrm{sc}^{([1],2)}$ is a $\Sigma^{\boldsymbol{\mathcal{A}}}$-homomorphism, according to Proposition~\ref{PDUCatHom}; the fourth equality follows from the fact that, by induction, $
\mathrm{sc}^{([1],2)}_{t}(
\mathrm{ip}^{(2,X)@}_{t}(
T^{(2)'}(M)))
=
\mathrm{sc}^{([1],2)}_{t}(
\mathrm{ip}^{(2,X)@}_{t}(
T^{(2)'}(N)))$; the fifth equality follows from the fact that  $\mathrm{sc}^{([1],2)}$ is a $\Sigma^{\boldsymbol{\mathcal{A}}}$-homomorphism, according to Proposition~\ref{PDUCatHom};  the sixth equality follows from the fact that $\mathrm{ip}^{(2,X)@}$ is a $\Sigma^{\boldsymbol{\mathcal{A}}^{(2)}}$-homomorphism, according to Definition~\ref{DDIp}; finally, the last equivalence follows from the description of the second-order elementary translation $T^{(2)}$.

The case of the $([1],2)$-target follows by a similar argument.

The remaining case, i.e., the case for which $s=t$, there is a second-order path $\mathfrak{P}^{(2)}\in\mathrm{Pth}_{\boldsymbol{\mathcal{A}}^{(2)},s}$ and a second-order translation $T^{(2)'}$ in $\mathrm{Tl}_{t}(\mathrm{T}_{\Sigma^{\boldsymbol{\mathcal{A}}^{(2)}}}(X))_{s}$ of height $m$  such that  such that
\[
T^{(2)}=\mathrm{CH}^{(2)}_{s}\left(\mathfrak{P}^{(2)}\right)
\circ^{0\mathbf{T}_{\Sigma^{\boldsymbol{\mathcal{A}}^{(2)}}}(X)}_{s}
T^{(2)'}
\]
follows by a similar argument to the proof presented above.

Now, regarding the other options for the second-order  translation $T^{(2)}$, it could be the case that $s=t$ and there is a second-order path $\mathfrak{P}^{(2)}\in\mathrm{Pth}_{\boldsymbol{\mathcal{A}}^{(2)},s}$ and a second-order translation $T^{(2)'}$ in $\mathrm{Tl}_{t}(\mathrm{T}_{\Sigma^{\boldsymbol{\mathcal{A}}^{(2)}}}(X))_{s}$ of height $m$  such that  such that
\[
T^{(2)}=T^{(2)'}\circ^{1\mathbf{T}_{\Sigma^{\boldsymbol{\mathcal{A}}^{(2)}}}(X)}_{s}\mathrm{CH}^{(2)}_{s}\left(\mathfrak{P}^{(2)}\right).
\]

For this case, the following chain of equivalences hold
\begin{flushleft}
$T^{(2)}(M)$ is a second-order path term
\allowdisplaybreaks
\begin{align*}
\Leftrightarrow\quad & T^{(2)'}(M)\circ^{1\mathbf{T}_{\Sigma^{\boldsymbol{\mathcal{A}}^{(2)}}}(X)}_{s}\mathrm{CH}^{(2)}_{s}\left(\mathfrak{P}^{(2)}\right)\mbox{ is a second-order path term}
\tag{1}
\\
\Leftrightarrow\quad &
\mathrm{ip}^{(2,X)@}_{s}\left(
T^{(2)'}(M)\circ^{1\mathbf{T}_{\Sigma^{\boldsymbol{\mathcal{A}}^{(2)}}}(X)}_{s}\mathrm{CH}^{(2)}_{s}\left(\mathfrak{P}^{(2)}\right)
\right)
\mbox{ is a second-order path}
\tag{2}
\\
\Leftrightarrow\quad &
\mathrm{ip}^{(2,X)@}_{s}\left(
T^{(2)'}(M)\right)\circ^{1\mathbf{Pth}_{\boldsymbol{\mathcal{A}}^{(2)}}}_{s}
\mathrm{ip}^{(2,X)@}_{s}\left(
\mathrm{CH}^{(2)}_{s}\left(\mathfrak{P}^{(2)}\right)
\right)
\\&\qquad\qquad\qquad\qquad\qquad\qquad\qquad\qquad\qquad\qquad\quad
\mbox{ is a second-order path}
\tag{3}
\\
\Leftrightarrow\quad &
\begin{cases}
\mathrm{ip}^{(2,X)@}_{s}\left(
T^{(2)'}(M)\right)\mbox{ is a secod-order path}
\\
\mathrm{sc}^{([1],2)}_{s}\left(
\mathrm{ip}^{(2,X)@}_{s}\left(
T^{(2)'}(M)\right)\right)
=
\mathrm{tg}^{([1],2)}_{s}\left(
\mathrm{ip}^{(2,X)@}_{s}\left(
\mathrm{CH}^{(2)}_{s}\left(\mathfrak{P}^{(2)}\right)
\right)\right)
\end{cases}
\tag{4}
\\
\Leftrightarrow\quad &
\begin{cases}
T^{(2)'}(M)\mbox{ is a second-order path term}
\\
\mathrm{sc}^{([1],2)}_{s}\left(
\mathrm{ip}^{(2,X)@}_{s}\left(
T^{(2)'}(M)\right)\right)
=
\mathrm{tg}^{([1],2)}_{s}\left(
\mathrm{ip}^{(2,X)@}_{s}\left(
\mathrm{CH}^{(2)}_{s}\left(\mathfrak{P}^{(2)}\right)
\right)\right)
\end{cases}
\tag{5}
\\
\Leftrightarrow\quad &
\begin{cases}
T^{(2)'}(N)\mbox{ is a second-order path term}
\\
\mathrm{sc}^{([1],2)}_{s}\left(
\mathrm{ip}^{(2,X)@}_{s}\left(
T^{(2)'}(N)\right)\right)
=
\mathrm{tg}^{([1],2)}_{s}\left(
\mathrm{ip}^{(2,X)@}_{s}\left(
\mathrm{CH}^{(2)}_{s}\left(\mathfrak{P}^{(2)}\right)
\right)\right)
\end{cases}
\tag{6}
\\
\Leftrightarrow\quad &
\begin{cases}
\mathrm{ip}^{(2,X)@}_{s}\left(
T^{(2)'}(N)\right)\mbox{ is a second-order path}
\\
\mathrm{sc}^{([1],2)}_{s}\left(
\mathrm{ip}^{(2,X)@}_{s}\left(
T^{(2)'}(N)\right)\right)
=
\mathrm{tg}^{([1],2)}_{s}\left(
\mathrm{ip}^{(2,X)@}_{s}\left(
\mathrm{CH}^{(2)}_{s}\left(\mathfrak{P}^{(2)}\right)
\right)\right)
\end{cases}
\tag{7}
\\
\Leftrightarrow\quad &
\mathrm{ip}^{(2,X)@}_{s}\left(
T^{(2)'}(N)\right)\circ^{1\mathbf{Pth}_{\boldsymbol{\mathcal{A}}^{(2)}}}_{s}
\mathrm{ip}^{(2,X)@}_{s}\left(
\mathrm{CH}^{(2)}_{s}\left(\mathfrak{P}^{(2)}\right)
\right)
\\&\qquad\qquad\qquad\qquad\qquad\qquad\qquad\qquad\qquad\qquad\quad
\mbox{ is a second-order path}
\tag{8}
\\
\Leftrightarrow\quad &
\mathrm{ip}^{(2,X)@}_{s}\left(
T^{(2)'}(N)\circ^{1\mathbf{T}_{\Sigma^{\boldsymbol{\mathcal{A}}^{(2)}}}(X)}_{s}\mathrm{CH}^{(2)}_{s}\left(\mathfrak{P}^{(2)}\right)
\right)
\mbox{ is a second-order path}
\tag{9}
\\
\Leftrightarrow\quad & T^{(2)'}(N)\circ^{1\mathbf{T}_{\Sigma^{\boldsymbol{\mathcal{A}}^{(2)}}}(X)}_{s}\mathrm{CH}^{(2)}_{s}\left(\mathfrak{P}^{(2)}\right)\mbox{ is a second-order path term}
\tag{10}
\\
\Leftrightarrow\quad &
T^{(2)}(N) \mbox{ is a second-order path term}.
\tag{11}
\end{align*}
\end{flushleft}

In the just stated chain of equivalences, the first equivalence follows from the description of the second-order  translation $T^{(2)}$; the second equivalence follows from Proposition~\ref{PDPT}; 	the third equivalence follows from the fact that $\mathrm{ip}^{(2,X)@}$ is a $\Sigma^{\boldsymbol{\mathcal{A}}^{(2)}}$-homomorphism, according to Definition~\ref{DDIp}; the fourth equivalence follows from the definition of the $1$-composition operation in $\mathbf{Pth}_{\boldsymbol{\mathcal{A}}^{(2)}}$ according to Proposition~\ref{PDPthDCatAlg}; 
the fifth equivalence follows from Proposition~\ref{PDPT};
the  sixth equivalence follows by induction. Moreover,  the equality $
\mathrm{sc}^{([1],2)}_{t}(
\mathrm{ip}^{(2,X)@}_{t}(
T^{(2)'}(M)))
=
\mathrm{sc}^{([1],2)}_{t}(
\mathrm{ip}^{(2,X)@}_{t}(
T^{(2)'}(N)))$ also holds by induction; the seventh equivalence follows from Proposition~\ref{PDPT}; the eighth equivalence follows from the definition of the $1$-composition operation in $\mathbf{Pth}_{\boldsymbol{\mathcal{A}}^{(2)}}$ according to Proposition~\ref{PDPthDCatAlg}; the ninth equivalence follows from the fact that $\mathrm{ip}^{(2,X)@}$ is a $\Sigma^{\boldsymbol{\mathcal{A}}^{(2)}}$-homomorphism, according to Definition~\ref{DDIp}; the tenth equivalence follows from Proposition~\ref{PDPT}; finally, the last equivalence follows from the description of the second-order translation $T^{(2)}$.

Assume that either $T^{(2)}(M)$ or $T^{(2)}(N)$ is a second-order path term.
Regarding the $([1],2)$-source, the following chain of equalities holds
\begin{flushleft}
$\mathrm{sc}^{([1],2)}_{s}\left(\mathrm{ip}^{(2,X)@}_{s}\left(
T^{(2)}(M)\right)\right)$
\allowdisplaybreaks
\begin{align*}
&=\mathrm{sc}^{([1],2)}_{s}\left(\mathrm{ip}^{(2,X)@}_{s}\left(
T^{(2)'}(M)\circ^{1\mathbf{T}_{\Sigma^{\boldsymbol{\mathcal{A}}^{(2)}}}(X)}_{s}\mathrm{CH}^{(2)}_{s}\left(\mathfrak{P}^{(2)}\right)\right)\right)
\tag{1}
\\&=
\mathrm{sc}^{([1],2)}_{s}\left(
\mathrm{ip}^{(2,X)@}_{s}\left(
T^{(2)'}(M)\right)
\circ^{1\mathbf{Pth}_{\boldsymbol{\mathcal{A}}^{(2)}}}_{s}
\mathrm{ip}^{(2,X)@}_{s}\left(
\mathrm{CH}^{(2)}_{s}\left(
\mathfrak{P}^{(2)}\right)
\right)\right)
\tag{2}
\\&=
\mathrm{sc}^{([1],2)}_{s}\left(
\mathrm{ip}^{(2,X)@}_{s}\left(
\mathrm{CH}^{(2)}_{s}\left(
\mathfrak{P}^{(2)}\right)
\right)\right)
\tag{3}
\\&=
\mathrm{sc}^{([1],2)}_{s}\left(
\mathrm{ip}^{(2,X)@}_{s}\left(
T^{(2)'}(N)\right)
\circ^{1\mathbf{Pth}_{\boldsymbol{\mathcal{A}}^{(2)}}}_{s}
\mathrm{ip}^{(2,X)@}_{s}\left(
\mathrm{CH}^{(2)}_{s}\left(
\mathfrak{P}^{2}\right)
\right)\right)
\tag{4}
\\&=\mathrm{sc}^{([1],2)}_{s}\left(\mathrm{ip}^{(2,X)@}_{s}\left(
T^{(2)'}(N)\circ^{1\mathbf{T}_{\Sigma^{\boldsymbol{\mathcal{A}}^{(2)}}}(X)}_{s}\mathrm{CH}^{(2)}_{s}\left(\mathfrak{P}^{(2)}\right)\right)\right)
\tag{5}
\\&=
\mathrm{sc}^{([1],2)}_{s}\left(\mathrm{ip}^{(2,X)@}_{s}\left(
T^{(2)}(N)\right)\right).
\tag{6}
\end{align*}
\end{flushleft}

In the just stated chain of equalities, the first equality follows from the description of the second-order  translation $T^{(2)}$; the second equality follows from the fact that $\mathrm{ip}^{(2,X)@}$ is a $\Sigma^{\boldsymbol{\mathcal{A}}^{(2)}}$-homomorphism, according to Definition~\ref{DDIp}; the third equality follows from the definition of the $1$-composition operation in $\mathbf{Pth}_{\boldsymbol{\mathcal{A}}^{(2)}}$ according to Proposition~\ref{PDPthDCatAlg}; the fourth equality follows from the definition of the $1$-composition operation in $\mathbf{Pth}_{\boldsymbol{\mathcal{A}}^{(2)}}$ according to Proposition~\ref{PDPthDCatAlg}; the fifth equality follows from the fact that $\mathrm{ip}^{(2,X)@}$ is a $\Sigma^{\boldsymbol{\mathcal{A}}^{(2)}}$-homomorphism, according to Definition~\ref{DDIp}; finally, the last equivalence follows from the description of the second-order elementary translation $T^{(2)}$.

Now, regarding the $([1],2)$-target, the following chain of equalities holds
\begin{flushleft}
$\mathrm{tg}^{([1],2)}_{s}\left(\mathrm{ip}^{(2,X)@}_{s}\left(
T^{(2)}(M)\right)\right)$
\allowdisplaybreaks
\begin{align*}
&=\mathrm{tg}^{([1],2)}_{s}\left(\mathrm{ip}^{(2,X)@}_{s}\left(
T^{(2)'}(M)\circ^{1\mathbf{T}_{\Sigma^{\boldsymbol{\mathcal{A}}^{(2)}}}(X)}_{s}\mathrm{CH}^{(2)}_{s}\left(\mathfrak{P}^{(2)}\right)\right)\right)
\tag{1}
\\&=
\mathrm{tg}^{([1],2)}_{s}\left(
\mathrm{ip}^{(2,X)@}_{s}\left(
T^{(2)'}(M)\right)
\circ^{1\mathbf{Pth}_{\boldsymbol{\mathcal{A}}^{(2)}}}_{s}
\mathrm{ip}^{(2,X)@}_{s}\left(
\mathrm{CH}^{(2)}_{s}\left(
\mathfrak{P}^{(2)}\right)
\right)\right)
\tag{2}
\\&=
\mathrm{tg}^{([1],2)}_{s}\left(
\mathrm{ip}^{(2,X)@}_{s}\left(
T^{(2)'}(M)\right)
\right)
\tag{3}
\\&=
\mathrm{tg}^{([1],2)}_{s}\left(
\mathrm{ip}^{(2,X)@}_{s}\left(
T^{(2)'}(N)\right)
\right)
\tag{4}
\\&=
\mathrm{tg}^{([1],2)}_{s}\left(
\mathrm{ip}^{(2,X)@}_{s}\left(
T^{(2)'}(N)\right)
\circ^{1\mathbf{Pth}_{\boldsymbol{\mathcal{A}}^{(2)}}}_{s}
\mathrm{ip}^{(2,X)@}_{s}\left(
\mathrm{CH}^{(2)}_{s}\left(
\mathfrak{P}^{(2)}\right)
\right)\right)
\tag{5}
\\&=\mathrm{tg}^{([1],2)}_{s}\left(\mathrm{ip}^{(2,X)@}_{s}\left(
T^{(2)'}(N)\circ^{1\mathbf{T}_{\Sigma^{\boldsymbol{\mathcal{A}}^{(2)}}}(X)}_{s}\mathrm{CH}^{(2)}_{s}\left(\mathfrak{P}^{(2)}\right)\right)\right)
\tag{6}
\\&=
\mathrm{tg}^{([1],2)}_{s}\left(\mathrm{ip}^{(2,X)@}_{s}\left(
T^{(2)}(N)\right)\right).
\tag{7}
\end{align*}
\end{flushleft}

In the just stated chain of equalities, the first equality follows from the description of the second-order  translation $T^{(2)}$; the second equality follows from the fact that $\mathrm{ip}^{(2,X)@}$ is a $\Sigma^{\boldsymbol{\mathcal{A}}^{(2)}}$-homomorphism, according to Definition~\ref{DDIp}; the third equality follows from the definition of the $1$-composition operation in $\mathbf{Pth}_{\boldsymbol{\mathcal{A}}^{(2)}}$ according to Proposition~\ref{PDPthDCatAlg}; the fourth
equality follows by induction. Note that, according to the above proof,  since $T^{(2)}(M)$ is a path term, then $T^{(2)'}(M)$ is a path term as well, which implies, by induction, that  $T^{(2)'}(N)$ is a path term. The implication also holds in the other direction. In particular,  $
\mathrm{tg}^{([1],2)}_{t}(
\mathrm{ip}^{(2,X)@}_{t}(
T^{(2)'}(M)))
=
\mathrm{tg}^{([1],2)}_{t}(
\mathrm{ip}^{(2,X)@}_{t}(
T^{(2)'}(N)))$; the fifth equality follows from the definition of the $1$-composition operation in $\mathbf{Pth}_{\boldsymbol{\mathcal{A}}^{(2)}}$ according to Proposition~\ref{PDPthDCatAlg}; the sixth equality follows from the fact that $\mathrm{ip}^{(2,X)@}$ is a $\Sigma^{\boldsymbol{\mathcal{A}}^{(2)}}$-homomorphism, according to Definition~\ref{DDIp}; finally, the last equivalence follows from the description of the second-order  translation $T^{(2)}$.

The remaining case, i.e., the case for which  $s=t$ and there is a second-order path $\mathfrak{P}^{(2)}\in\mathrm{Pth}_{\boldsymbol{\mathcal{A}}^{(2)},s}$ and a second-order translation $T^{(2)'}$ in $\mathrm{Tl}_{t}(\mathrm{T}_{\Sigma^{\boldsymbol{\mathcal{A}}^{(2)}}}(X))_{s}$ of height $m$  such that  such that
\[
T^{(2)}=\mathrm{CH}^{(2)}_{s}\left(\mathfrak{P}^{(2)}\right)\circ^{1\mathbf{T}_{\Sigma^{\boldsymbol{\mathcal{A}}^{(2)}}}(X)}_{s} T^{(2)'},
\]
follows by a similar argument to the proof presented above.

This complete the proof.
\end{proof}

\begin{corollary}
\label{CDTransWD} 
Let $s, t$ be sorts in $S$, $T^{(2)}$ a second-order translation in 
$\mathrm{Tl}_{t}(\mathrm{T}_{\Sigma^{\boldsymbol{\mathcal{A}}^{(2)}}}(X))_{s}$ and $M$ and $N$ be second-order path terms in
$\mathrm{PT}_{\boldsymbol{\mathcal{A}}^{(2)},t}$ satisfying that 
 $$
\left(
\mathrm{ip}^{(2,X)@}_{t}\left(
M
\right), 
\mathrm{ip}^{(2,X)@}_{t}\left(
N
\right)\right)
\in\mathrm{Ker}\left(
\mathrm{sc}^{([1],2)}
\right)_{t}
\cap\mathrm{Ker}\left(
\mathrm{tg}^{([1],2)}
\right)_{t}.
$$

If either $T^{(2)}(M)$ or $T^{(2)}(N)$ is a second-order path term in $\mathrm{PT}_{\boldsymbol{\mathcal{A}}^{(2)},s}$, then 
\[
\left(
\mathrm{ip}^{(2,X)@}_{s}\left(
T^{(2)}(M)\right),
\mathrm{ip}^{(2,X)@}_{s}\left(
T^{(2)}(N)\right)
\right)
\in\mathrm{Ker}\left(
\mathrm{sc}^{(0,2)}
\right)_{s}
\cap\mathrm{Ker}\left(
\mathrm{tg}^{(0,2)}
\right)_{s}.
\]
\end{corollary}
\begin{proof}
It follows from Lemma~\ref{LDTransWD} and Definition~\ref{DDScTgZ}.
\end{proof}

\begin{restatable}{proposition}{PDTransWD}
\label{PDTransWD} Let $s, t$ be sorts in $S$, $T^{(2)}$ a second-order translation in 
$\mathrm{Tl}_{t}(\mathrm{T}_{\Sigma^{\boldsymbol{\mathcal{A}}^{(2)}}}(X))_{s}$ and  $P$ and $P'$ be second-order path terms in
$\mathrm{PT}_{\boldsymbol{\mathcal{A}}^{(2)},t}$ satisfying that $(P,P')\in \Theta^{\llbracket 2\rrbracket}_{t}$.

If either $T^{(2)}(P)$ or $T^{(2)}(P')$ is a second-order path term in $\mathrm{PT}_{\boldsymbol{\mathcal{A}}^{(2)},s}$, then 
\[
\left(
T^{(2)}(P),T^{(2)}(P')
\right)
\in\Theta^{\llbracket 2\rrbracket}_{s}.
\]
\end{restatable}
\begin{proof}
It follows from the fact that $\Theta^{\llbracket 2\rrbracket}$ is a $\Sigma^{\boldsymbol{\mathcal{A}}^{(2)}}$-congruence on $\mathbf{T}_{\Sigma^{\boldsymbol{\mathcal{A}}^{(2)}}}$, according to Definition~\ref{DDVCong}, and Proposition~\ref{PTransCong}.
\end{proof}


\part{Morphisms of rewriting systems}
\chapter{Morphisms of first-order many-sorted rewriting systems}\label{S3A}

Let us recall that a first-order many-sorted rewriting system is an ordered pair $\boldsymbol{\mathcal{A}}^{(1)}=(\boldsymbol{\mathcal{A}}^{(0)}, \mathcal{A}^{(1)})$, where $\boldsymbol{\mathcal{A}}^{(0)} = (S, \Sigma, X)$ is a zeroth-order $S$-sorted rewriting system, i.e., $S$ is a set of sorts, $\Sigma$ is a $S$-sorted signature and $X$ an $S$-sorted set, and $\mathcal{A}^{(1)}$ is a subset of $\mathrm{Rwr}(\mathbf{\Sigma}, X) = (\T_{\Sigma}(X)_{s}^{2})_{s \in S}$, where $\T_{\Sigma}(X)$ is the underlying $A$-sorted set of $\mathbf{T}_{\Sigma}(X)$, the free many-sorted $\Sigma$-algebra on $X$. For further details, see Definition~\ref{DRewSys}.

In this chapter, given two first-order many-sorted rewriting systems $\boldsymbol{\mathcal{A}}^{(1)}=(\boldsymbol{\mathcal{A}}^{(0)}, \mathcal{A}^{(1)})$ and  $\boldsymbol{\mathcal{B}}^{(1)}=(\boldsymbol{\mathcal{B}}^{(0)}, \mathcal{B}^{(1)})$, we define the notion of first-order morphism $\mathbf{f}^{(1)}$ from  $\boldsymbol{\mathcal{A}}^{(1)}$ to  $\boldsymbol{\mathcal{B}}^{(1)}$, where $\mathbf{f}^{(1)}=(\mathbf{f}^{(0)},f^{(1)})$, $\mathbf{f}^{(0)}=(\varphi,c,f^{(0)})$ is a zeroth-order morphism from $\boldsymbol{\mathcal{A}}^{(0)}=(S,\Sigma,X)$ to $\boldsymbol{\mathcal{B}}^{(0)}=(T,\Lambda,Y)$, see Definition~\ref{DRws0Mor}, and $f^{(1)}$ is an $S$-sorted mapping from $\mathcal{A}^{(1)}$ to $\mathrm{Pth}_{\boldsymbol{\mathcal{B}}^{(1)},\varphi}$ with some compatibility properties with sources and targets. We next define a structure of $\Sigma$-algebra on $\mathrm{Pth}_{\boldsymbol{\mathcal{B}}^{(1)}, \varphi}$ and $[ \mathrm{Pth}_{\boldsymbol{\mathcal{B}}^{(1)}, \varphi} ]$, which we will denote by $\mathbf{Pth}_{\boldsymbol{\mathcal{B}}^{(1)}}^{\mathbf{f}^{(1)}(0,1)}$ and $[ \mathbf{Pth}_{\boldsymbol{\mathcal{B}}^{(1)}}^{\mathbf{f}^{(1)}(0,1)} ]$, respectively. Then, given a first-order morphism $\mathbf{f}^{(1)}$ from $\boldsymbol{\mathcal{A}}^{(1)}$ to $\boldsymbol{\mathcal{B}}^{(1)}$, we consider an extension of the mapping $f^{(1)}$ to the $S$-sorted set $\mathrm{Pth}_{\boldsymbol{\mathcal{A}}^{(1)}}$, which we have called the path extension mapping, and denote it by $f^{(1)\flat}$.  Finally, we prove that the path extension mapping is a $\Sigma$-homomorphism.

We first introduce the notion of morphism from a first-order many-sorted rewriting to another.

\begin{definition}\label{DRws1Mor}
A \emph{morphism} of first-order many-sorted rewriting systems, or simply a \emph{first-order morphism}, from $\boldsymbol{\mathcal{A}}^{(1)}=(\boldsymbol{\mathcal{A}}^{(0)}, \mathcal{A}^{(1)})$ to $\boldsymbol{\mathcal{B}}^{(1)}=(\boldsymbol{\mathcal{B}}^{(0)}, \mathcal{B}^{(1)})$, is an ordered triple $(\boldsymbol{\mathcal{A}}^{(1)}, \mathbf{f}^{(1)}, \boldsymbol{\mathcal{B}}^{(1)})$, denoted by $\mathbf{f}^{(1)} \colon \boldsymbol{\mathcal{A}}^{(1)} \mor \boldsymbol{\mathcal{B}}^{(1)}$ for short, in which $\mathbf{f}^{(1)}=(\varphi, c, (f^{(i)})_{i\in 2})$ is the ordered pair where 
\begin{enumerate}
\item
$\mathbf{f}^{(0)} = (\varphi, c, f^{(0)})$, the \emph{underlying zeroth-order morphism} of $\mathbf{f}^{(1)}$, is a zeroth-order morphism from $\boldsymbol{\mathcal{A}}^{(0)}=(S,\Sigma,X)$ to $\boldsymbol{\mathcal{B}}^{(0)}=(T,\Lambda,Y)$, as introduced in Definition~\ref{DRws0Mor},
$$
\mathbf{f}^{(0)} \colon
\boldsymbol{\mathcal{A}}^{(0)}
\mor
\boldsymbol{\mathcal{B}}^{(0)}
$$
\item 
$f^{(1)} \colon \mathcal{A}^{(1)} \mor \mathrm{Pth}_{\boldsymbol{\mathcal{B}}^{(1)},\varphi}$ is an $S$-sorted mapping, where $\mathrm{Pth}_{\boldsymbol{\mathcal{B}}^{(1)},\varphi}$ is the $S$-sorted set $(\mathrm{Pth}_{\boldsymbol{\mathcal{B}}^{(1)},\varphi(s)})_{s\in S}$, satisfying that, for every $s\in S$ and every rewrite rule $\mathfrak{p}=(M,N)\in\mathcal{A}^{(1)}_{s}$, we have that
$$
f_{s}^{(1)}\left(\mathfrak{p}\right)
\in
\mathrm{Pth}_{\boldsymbol{\mathcal{B}}^{(1)}, \varphi(s)}\left(
f^{(0)\sharp}_{s}\left(
M
\right), 
f^{(0)\sharp}_{s}\left(
N
\right)
\right),
$$
or, what is equivalent,
\begin{align*}
\mathrm{sc}^{(0,1)}_{\boldsymbol{\mathcal{B}}^{(1)},\varphi(s)}\left(
f_{s}^{(1)}\left(
\mathfrak{p}
\right)
\right)
&=
f^{(0)\sharp}_{s}\left(
M
\right)
&&\mbox{ and}&
\mathrm{tg}^{(0,1)}_{\boldsymbol{\mathcal{B}}^{(1)},\varphi(s)}\left(
f_{s}^{(1)}\left(
\mathfrak{p}
\right)
\right)
&=
f^{(0)\sharp}_{s}\left(
N
\right).
\end{align*}
\end{enumerate}

The alternative notation $\mathbf{f}^{(1)} = (\mathbf{f}^{(0)}, f^{(1)})$ will also be used.

\end{definition}

The following remark follows directly from Proposition~\ref{PFunSig}.

\begin{remark}\label{RDSigmaAlg}
Let $\mathbf{f}^{(1)}=(\varphi, c, (f^{(i)})_{i\in 2})$ be a first-order morphism from $\boldsymbol{\mathcal{A}}^{(1)}$ to $\boldsymbol{\mathcal{B}}^{(1)}$.
Then the $S$-sorted sets $\mathrm{Pth}_{\boldsymbol{\mathcal{B}}^{(1)}, \varphi}$ and $\left[\mathrm{Pth}_{\boldsymbol{\mathcal{B}}^{(1)}}\right]_{\varphi}$ are equipped, in a natural way, with a structure of $\Sigma$-algebra, namely, $\mathbf{c}_{\mathfrak{d}}^{\ast}(\mathbf{Pth}_{\boldsymbol{\mathcal{B}}^{(1)}}^{(0,1)})$ and $\mathbf{c}_{\mathfrak{d}}^{\ast}([\mathbf{Pth}_{\boldsymbol{\mathcal{B}}^{(1)}}^{(0,1)}])$.
In order to unify the notation, we will denote these $\Sigma$-algebras by $\mathbf{Pth}_{\boldsymbol{\mathcal{B}}^{(1)}}^{\mathbf{f}^{(1)}(0,1)}$ and $[\mathbf{Pth}_{\boldsymbol{\mathcal{B}}^{(1)}}^{\mathbf{f}^{(1)}(0,1)}]$, respectively.
Moreover, the $S$-sorted mapping $\mathrm{pr}_{\boldsymbol{\mathcal{B}}^{(1)}, \varphi}^{\mathrm{Ker}(\mathrm{CH}^{(1)})} = \mathbf{c}_{\mathfrak{d}}^{\ast}(\mathrm{pr}_{\boldsymbol{\mathcal{B}}^{(1)}}^{\mathrm{Ker}(\mathrm{CH}^{(1)})})$ is a $\Sigma$-homomorphism from  $\mathbf{Pth}_{\boldsymbol{\mathcal{B}}^{(1)}}^{\mathbf{f}^{(1)}(0,1)}$ to $[\mathbf{Pth}_{\boldsymbol{\mathcal{B}}^{(1)}}^{\mathbf{f}^{(1)}(0,1)}]$ by Propositions~\ref{PCHHom} and \ref{PFunSig}.
Similarly, $\mathrm{sc}^{(0,1)}_{\boldsymbol{\mathcal{B}}^{(1)}, \varphi} = \mathbf{c}_{\mathfrak{d}}^{\ast}(\mathrm{sc}^{(0,1)}_{\boldsymbol{\mathcal{B}}^{(1)}})$ and $\mathrm{tg}^{(0,1)}_{\boldsymbol{\mathcal{B}}^{(1)}, \varphi} = \mathbf{c}_{\mathfrak{d}}^{\ast}(\mathrm{tg}^{(0,1)}_{\boldsymbol{\mathcal{B}}^{(1)}})$ are $\Sigma$-homomorphisms from $\mathbf{Pth}_{\boldsymbol{\mathcal{B}}^{(1)}}^{\mathbf{f}^{(1)}(0,1)}$ to $\mathbf{c}_{\mathfrak{d}}^{\ast}(\mathbf{T}_{\Lambda}(Y))$ by Propositions~\ref{PHom} and \ref{PFunSig}.
Finally, $\mathrm{ip}^{(1,0)\sharp}_{\boldsymbol{\mathcal{B}}^{(1)}, \varphi} = \mathbf{c}_{\mathfrak{d}}^{\ast} (\mathrm{ip}^{(1,0)\sharp}_{\boldsymbol{\mathcal{B}}^{(1)}})$ is a  $\Sigma$-homomorphism from $\mathbf{c}_{\mathfrak{d}}^{\ast}(\mathbf{T}_{\Lambda}(Y))$ to $\mathbf{Pth}_{\boldsymbol{\mathcal{B}}^{(1)}}^{\mathbf{f}^{(1)}(0,1)}$ by Propositions~\ref{PIpHom} and \ref{PFunSig}.
\end{remark}

\section{Path extension mapping of a morphism}

In this section we prove that, for a first-order morphism $\mathbf{f}^{(1)}$, there exists an extension of $f^{(1)}$, the \emph{path extension mapping}, to the set of paths $\mathrm{Pth}_{\boldsymbol{\mathcal{A}}^{(1)}}$ and show its uniqueness for some conditions.

\begin{proposition}
\label{PPthExt}
Let $\mathbf{f}^{(1)}=(\varphi, c, (f^{(i)})_{i\in 2})$ be a first-order morphism from $\boldsymbol{\mathcal{A}}^{(1)}$ to $\boldsymbol{\mathcal{B}}^{(1)}$. Then there exists an $S$-sorted mapping $f^{(1)\flat}$ from $\mathrm{Pth}_{\boldsymbol{\mathcal{A}}^{(1)}}$ to $\mathrm{Pth}_{\boldsymbol{\mathcal{B}}^{(1)}, \varphi}$, which we call the \emph{path extension mapping of $f^{(1)}$}, satisfying that
\begin{multicols}{2}
\begin{enumerate}
\item $\mathrm{sc}^{(0,1)}_{\boldsymbol{\mathcal{B}}^{(1)},\varphi}\circ f^{(1)\flat}=f^{(0)\sharp}\circ \mathrm{sc}^{(0,1)}_{\boldsymbol{\mathcal{A}}^{(1)}}$;
\item $\mathrm{tg}^{(0,1)}_{\boldsymbol{\mathcal{B}}^{(1)},\varphi}\circ f^{(1)\flat}=f^{(0)\sharp}\circ \mathrm{tg}^{(0,1)}_{\boldsymbol{\mathcal{A}}^{(1)}}$;
\item
$f^{(1)\flat} \circ \mathrm{ip}_{\boldsymbol{\mathcal{A}}^{(1)}}^{(1, 0)\sharp} = \mathrm{ip}_{\boldsymbol{\mathcal{B}}^{(1)}, \varphi}^{(1, 0)\sharp} \circ f^{(0)\sharp}$;
\item
$f^{(1)\flat} \circ \mathrm{ech}^{(1,\mathcal{A}^{(1)})}_{\boldsymbol{\mathcal{A}}^{(1)}} = f^{(1)}$.
\end{enumerate}
\end{multicols}
\end{proposition}

\begin{proof}
Let us define the $S$-sorted mapping $f^{(1)\flat}$ by Artinian recursion on $(\coprod \mathrm{Pth}_{\boldsymbol{\mathcal{A}}^{(1)}}, \leq_{\mathbf{Pth}_{\boldsymbol{\mathcal{A}}^{(1)}}})$ as follows.

{\sffamily Base step of the Artinian recursion.}

Let $(\mathfrak{P}, s)$ be a minimal element of $(\coprod \mathrm{Pth}_{\boldsymbol{\mathcal{A}}^{(1)}}, \leq_{\mathbf{Pth}_{\boldsymbol{\mathcal{A}}^{(1)}}})$. Then, by Proposition~\ref{PMinimal}, the path $\mathfrak{P}$ is either~(1) an $(1,0)$-identity path or~(2) an echelon.

If~(1), i.e., if $\mathfrak{P}$ is an $(1,0)$-identity path, then $\mathfrak{P}=\mathrm{ip}^{(1,0)\sharp}_{\boldsymbol{\mathcal{A}}^{(1)}, s}(P)$ for some term $P \in \T_{\Sigma}(X)_{s}$. We define $f^{(1)\flat}_{s}(\mathfrak{P})$ to be the $(1,0)$-identity  path at $f_{s}^{(0)\sharp}(P)$ which is a term in $\T_{\Lambda}(Y)_{\varphi(s)}$, i.e.,
$$
f^{(1)\flat}_{s}(\mathfrak{P}) = \mathrm{ip}^{(1,0)\sharp}_{\boldsymbol{\mathcal{B}}^{(1)}, \varphi(s)}\left(
f_{s}^{(0)\sharp}\left(
P
\right)
\right).
$$

If~(2), i.e., if $\mathfrak{P}$ is an echelon associated to a rewrite rule $\mathfrak{p}=(M, N)$, that is, if $\mathfrak{P}$ has the form
$$
\xymatrix@C=55pt{
\mathfrak{P}: M
\ar[r]^-{\text{\Small{$(\mathfrak{p},\mathrm{id}^{\mathrm{T}_{\Sigma}(X)_{s}})$}}}
&
N
},
$$
then we define $f_{s}^{(1)\flat}\left(\mathfrak{P}\right)$ to be the image of $\mathfrak{p}$, the unique rewrite rule appearing in $\mathfrak{P}$, under the mapping $f_{s}^{(1)}$, i.e.,
$$
f^{(1)\flat}_{s}\left(\mathfrak{P}\right)
=
f_{s}^{(1)}\left(\mathfrak{p}\right).
$$

This completes the base step of the Artinian recursion.

{\sffamily Inductive step of the Artinian recursion.}

Let $(\mathfrak{P},s)$ be a non-minimal element of $(\coprod\mathrm{Pth}_{\boldsymbol{\mathcal{A}}^{(1)}}, \leq_{\mathbf{Pth}_{\boldsymbol{\mathcal{A}}^{(1)}}})$. We can assume that $\mathfrak{P}$ is a not a $(1,0)$-identity path, since those paths already have an image for the mapping $f^{(1)\flat}$. Let us suppose that, for every sort $t\in S$ and every path $\mathfrak{Q} \in\mathrm{Pth}_{\boldsymbol{\mathcal{A}}^{(1)},t}$, if $(\mathfrak{Q},t)<_{\mathbf{Pth}_{\boldsymbol{\mathcal{A}}^{(1)}}}(\mathfrak{P},s)$, then the value of the mapping $f^{(1)\flat}$ at $\mathfrak{Q}$, i.e., $f^{(1)\flat}_{t}(\mathfrak{Q})$, has already been defined. Moreover, assume that, for every $t \in S$ and every path $\mathfrak{Q}\in\mathrm{Pth}_{\boldsymbol{\mathcal{A}}^{(1)}, t}$, the definition of $f^{(1)\flat}_{t}(\mathfrak{Q})$ satisfies the following equalities
\allowdisplaybreaks
\begin{align*}
\mathrm{sc}_{\boldsymbol{\mathcal{B}}^{(1)},\varphi(t)}^{(0,1)}\left(f^{(1)\flat}_{t}\left(\mathfrak{Q}\right)\right) 
&=
f^{(0)\sharp}_{t}\left(\mathrm{sc}_{\boldsymbol{\mathcal{A}}^{(1)},t}^{(0,1)}\left(\mathfrak{Q}\right)\right);
\\
\mathrm{tg}_{\boldsymbol{\mathcal{B}}^{(1)},\varphi(t)}^{(0,1)}\left(f^{(1)\flat}_{t}\left(\mathfrak{Q}\right)\right)
&=
f^{(0)\sharp}_{t}\left(\mathrm{tg}_{\boldsymbol{\mathcal{A}}^{(1)},t}^{(0,1)}\left(\mathfrak{Q}\right)\right).
\end{align*}

By Lemma~\ref{LOrdI}, we have that $\mathfrak{P}$ is either~(1) a path of length strictly greater than one containing at least one echelon or~(2) an echelonless path.

If~(1), i.e., if $\mathfrak{P}$ is a path of length strictly greater than one containing at least one echelon, then let $i\in \bb{\mathfrak{P}}$ be the first index for which the one-step subpath $\mathfrak{P}^{i,i}$ of $\mathfrak{P}$ is an echelon. We consider different cases for $i$ according to the cases presented in Definition~\ref{DOrd}.

If $i=0$, we have that the pairs $(\mathfrak{P}^{0,0},s)$ and $(\mathfrak{P}^{1,\bb{\mathfrak{P}}-1},s)$ $\prec_{\mathbf{Pth}_{\boldsymbol{\mathcal{A}}^{(1)}}}$-precede the pair $(\mathfrak{P},s)$. Therefore, the values of the mapping $f^{(1)\flat}_{s}$ at $\mathfrak{P}^{0,0}$ and $\mathfrak{P}^{1,\bb{\mathfrak{P}}-1}$, respectively, have already been defined. In particular, the following chain of equalities holds
\allowdisplaybreaks
\begin{align*}
\mathrm{sc}^{(0,1)}_{\boldsymbol{\mathcal{B}}^{(1)},\varphi(s)}\left(
f_{s}^{(1)\flat}\left(
\mathfrak{P}^{1,\bb{\mathfrak{P}}-1}
\right)
\right)
&=
f^{(0)\sharp}_{s}\left(
\mathrm{sc}^{(0,1)}_{\boldsymbol{\mathcal{A}}^{(1)},s}\left(
\mathfrak{P}^{1,\bb{\mathfrak{P}}-1}
\right)
\right)
\tag{1}
\\
&=
f^{(0)\sharp}_{s}\left(
\mathrm{tg}^{(0,1)}_{\boldsymbol{\mathcal{A}}^{(1)},s}\left(
\mathfrak{P}^{0,0}
\right)
\right)
\tag{2}
\\
&=
\mathrm{tg}^{(0,1)}_{\boldsymbol{\mathcal{B}}^{(1)},\varphi(s)}\left(
f_{s}^{(1)\flat}\left(
\mathfrak{P}^{0,0}
\right)
\right).
\tag{3}
\end{align*}

The first equality follows from the assumption of the Artinian recursion;
the first equality follows from Proposition~\ref{PPthRecons};
finally, the last equality follows from the assumption of the Artinian recursion.

In this case, we set
\allowdisplaybreaks
\begin{align*}
f^{(1)\flat}_{s}\left(
\mathfrak{P}
\right)
&=
f^{(1)\flat}_{s}\left(
\mathfrak{P}^{1,\bb{\mathfrak{P}}-1}
\right)
\circ_{\varphi(s)}^{0\mathbf{Pth}_{\boldsymbol{\mathcal{B}}^{(1)}}}
f^{(1)\flat}_{s}\left(
\mathfrak{P}^{0,0}
\right).
\end{align*}
The above chain of equalities justifies that such composite exists.

If $i\neq 0$, we have that the pairs $(\mathfrak{P}^{0,i-1},s)$ and $(\mathfrak{P}^{i,\bb{\mathfrak{P}}-1},s)$ $\prec_{\mathbf{Pth}_{\boldsymbol{\mathcal{A}}^{(1)}}}$-precede the pair $(\mathfrak{P},s)$. Therefore, the values of the mapping $f^{(1)\flat}_{s}$ at $\mathfrak{P}^{0,i-1}$ and $\mathfrak{P}^{i,\bb{\mathfrak{P}}-1}$, respectively, have already been defined. In particular, the following chain of equalities holds
\allowdisplaybreaks
\begin{align*}
\mathrm{sc}^{(0,1)}_{\boldsymbol{\mathcal{B}}^{(1)},\varphi(s)}\left(
f_{s}^{(1)\flat}\left(
\mathfrak{P}^{i,\bb{\mathfrak{P}}-1}
\right)
\right)
&=
f^{(0)\sharp}_{s}\left(
\mathrm{sc}^{(0,1)}_{\boldsymbol{\mathcal{A}}^{(1)},s}\left(
\mathfrak{P}^{i,\bb{\mathfrak{P}}-1}
\right)
\right)
\tag{1}
\\
&=
f^{(0)\sharp}_{s}\left(
\mathrm{tg}^{(0,1)}_{\boldsymbol{\mathcal{A}}^{(1)},s}\left(
\mathfrak{P}^{0,i-1}
\right)
\right)
\tag{2}
\\
&=
\mathrm{tg}^{(0,1)}_{\boldsymbol{\mathcal{B}}^{(1)},\varphi(s)}\left(
f_{s}^{(1)\flat}\left(
\mathfrak{P}^{0,i-1}
\right)
\right).
\tag{3}
\end{align*}

The first equality follows from the assumption of the Artinian recursion; the first equality follows from Proposition~\ref{PPthRecons}; finally, the last equality follows from the assumption of the Artinian recursion.

In this case, we set
\allowdisplaybreaks
\begin{align*}
f^{(1)\flat}_{s}\left(
\mathfrak{P}
\right)
&=
f^{(1)\flat}_{s}\left(
\mathfrak{P}^{i,\bb{\mathfrak{P}}-1}
\right)
\circ_{\varphi(s)}^{0\mathbf{Pth}_{\boldsymbol{\mathcal{B}}^{(1)}}}
f^{(1)\flat}_{s}\left(
\mathfrak{P}^{0,i-1}
\right).
\end{align*}
The above chain of equalities justifies that such composite exists.

This finishes the definition of the value of the mapping $f^{(1)\flat}$ at a path of length strictly greater than one containing at least one echelon.

If~(2), i.e., if $\mathfrak{P}$ is an echelonless path in $\mathrm{Pth}_{\boldsymbol{\mathcal{A}}^{(1)},s}$, then the conditions for the path extraction algorithm, as stated in Lemma~\ref{LPthExtract}, are fulfilled. Then, by Lemma~\ref{LPthHeadCt}, there exists a unique word $\mathbf{s} \in S^{\star} - \{ \lambda \}$ and a unique operation symbol $\sigma \in \Sigma_{\mathbf{s}, s}$ associated to $\mathfrak{P}$. Let $(\mathfrak{P}_{j})_{j \in \bb{\mathbf{s}}}$ be the family of paths in $\mathrm{Pth}_{\boldsymbol{\mathcal{A}}^{(1)}, s}$ which, in virtue of Lemma~\ref{LPthExtract}, we can extract from $\mathfrak{P}$. Note that, for every $j \in \bb{\mathbf{s}}$, we have that $(\mathfrak{P}_{j}, s_{j}) \prec_{\mathbf{Pth}_{\boldsymbol{\mathcal{A}}^{(1)}}} (\mathfrak{P}, s)$. Therefore, for every $j \in \bb{\mathbf{s}}$, the value of the mapping $f^{(1)\flat}$ at $\mathfrak{P}_{j}$ has already been defined.

In this case, we set
$$
f^{(1)\flat}_{s}\left(
\mathfrak{P}
\right)
=
\sigma^{\mathbf{Pth}_{\boldsymbol{\mathcal{B}}^{(1)}}^{\mathbf{f}^{(1)}(0,1)}}\left(\left(
f^{(1)\flat}_{s_{j}}\left(
\mathfrak{P}_{j}
\right)
\right)_{j\in\bb{\mathbf{s}}}\right).
$$

This completes the definition of the mapping $f^{(1)\flat}$.

Now, we prove that the just defined mapping satisfies all the properties listed in the proposition.

\textsf{(1)} $\mathrm{sc}^{(0,1)}_{\boldsymbol{\mathcal{B}}^{(1)},\varphi}\circ f^{(1)\flat}=f^{(0)\sharp}\circ \mathrm{sc}^{(0,1)}_{\boldsymbol{\mathcal{A}}^{(1)}}$.

Let $s$ be a sort in $S$ and let $\mathfrak{P}$ be a path in $\mathrm{Pth}_{\boldsymbol{\mathcal{A}}^{(1)},s}$. We prove that
\begin{align*}
\mathrm{sc}_{\boldsymbol{\mathcal{B}}^{(1)},\varphi(s)}^{(0,1)}\left(
f_{s}^{(1)\flat}\left(
\mathfrak{P}
\right)
\right)
&=
f_{s}^{(0)\sharp}\left(
\mathrm{sc}_{\boldsymbol{\mathcal{A}}^{(1)},s}^{(0,1)}\left(
\mathfrak{P}
\right)
\right)
\end{align*}
by Artinian induction on $(\coprod\mathrm{Pth}_{\boldsymbol{\mathcal{A}}^{(1)}}, \leq_{\mathbf{Pth}_{\boldsymbol{\mathcal{A}}^{(1)}}})$.

{\sffamily Base step of the Artinian induction.}

Let $(\mathfrak{P}, s)$ be a minimal element of $(\coprod \mathrm{Pth}_{\boldsymbol{\mathcal{A}}^{(1)}}, \leq_{\mathbf{Pth}_{\boldsymbol{\mathcal{A}}^{(1)}}})$. Then, by Proposition~\ref{PMinimal}, the path $\mathfrak{P}$ is either~(1) an $(1,0)$-identity path or~(2) an echelon.

If~(1), i.e., if $\mathfrak{P}$ is an $(1,0)$-identity path, then $\mathfrak{P}=\mathrm{ip}^{(1,0)\sharp}_{\boldsymbol{\mathcal{A}}^{(1)}, s}(P)$ for some term $P \in \T_{\Sigma}(X)_{s}$.  Note that the following chain of equalities holds
\allowdisplaybreaks
\begin{align*}
\mathrm{sc}_{\boldsymbol{\mathcal{B}}^{(1)},\varphi(s)}^{(0,1)}\left(
f^{(1)\flat}_{s}\left(
\mathfrak{P}
\right)
\right) 
&= 
\mathrm{sc}_{\boldsymbol{\mathcal{B}}^{(1)},\varphi(s)}^{(0,1)}\left(
\mathrm{ip}_{\boldsymbol{\mathcal{B}}^{(1)},\varphi(s)}^{(1,0)\sharp}\left(
f_{s}^{(0)\sharp}\left(
P
\right)
\right)
\right)
\tag{1}
\\
&=
f_{s}^{(0)\sharp}\left(
P
\right)
\tag{2}
\\
&=
f_{s}^{(0)\sharp}\left(
\mathrm{sc}_{\boldsymbol{\mathcal{A}}^{(1)},s}^{(0,1)}\left(
\mathrm{ip}_{\boldsymbol{\mathcal{A}}^{(1)},s}^{(1,0)\sharp}\left(
P
\right)
\right)
\right)
\tag{3}
\\
&=
f^{(0)\sharp}_{s}\left(
\mathrm{sc}_{\boldsymbol{\mathcal{A}}^{(1)},s}^{(0,1)}
\left(
\mathfrak{P}
\right)
\right).
\tag{4}
\end{align*}

The first equality unravels the definition of $f^{(1)\flat}$;
the second and third equalities follow from the fact that, according to Proposition~\ref{PBasicEq}, $\mathrm{sc}_{\boldsymbol{\mathcal{B}}^{(1)},\varphi(s)}^{(0,1)} \circ \mathrm{ip}^{(1,0)\sharp}_{\boldsymbol{\mathcal{B}}^{(1)}, \varphi(s)}=\mathrm{id}^{\mathrm{T}_{\Lambda}(Y)_{\varphi(s)}}$ and $\mathrm{sc}_{\boldsymbol{\mathcal{A}}^{(1)},s}^{(0,1)} \circ \mathrm{ip}^{(1,0)\sharp}_{\boldsymbol{\mathcal{A}}^{(1)}, s} 
= 
\mathrm{id}^{\mathrm{T}_{\Sigma}(X)_{s}}$; 
finally, the fourth equality recovers the definition of $\mathfrak{P}$.

If~(2), i.e., if $\mathfrak{P}$ is an echelon associated to a rewrite rule $\mathfrak{p}=(M, N)$, the following chain of equalities holds
\allowdisplaybreaks
\begin{align*}
\mathrm{sc}_{\boldsymbol{\mathcal{B}}^{(1)},\varphi(s)}^{(0,1)}\left(
f^{(1)\flat}_{s}\left(
\mathfrak{P}
\right)
\right)
&=
\mathrm{sc}_{\boldsymbol{\mathcal{B}}^{(1)},\varphi(s)}^{(0,1)}\left(
f^{(1)}_{s}\left(
\mathfrak{p}
\right)
\right) 
\tag{1}
\\
&=
f^{(0)\sharp}_{s}\left(
M
\right)
\tag{2}
\\
&=
f^{(0)\sharp}_{s}\left(
\mathrm{sc}_{\boldsymbol{\mathcal{A}}^{(1)}, s}^{(0,1)} \left(
\mathfrak{P}
\right)
\right).
\tag{3}
\end{align*}

The first equality unravels the definition of $f^{(1)\flat}$;
the first equality follows from item~(2) of Definition~\ref{DRws1Mor};
the last equality follows from the definition of $\mathfrak{P}$.

This completes the base step of the Artinian induction.

{\sffamily Inductive step of the Artinian induction.}

Let $(\mathfrak{P},s)$ be a non-minimal element of $(\coprod\mathrm{Pth}_{\boldsymbol{\mathcal{A}}^{(1)}}, \leq_{\mathbf{Pth}_{\boldsymbol{\mathcal{A}}^{(1)}}})$. We can assume that $\mathfrak{P}$ is a not a $(1,0)$-identity path, since this case has already been proven. Let us suppose that, for every sort $t\in S$ and every path $\mathfrak{Q} \in\mathrm{Pth}_{\boldsymbol{\mathcal{A}}^{(1)},t}$, if $(\mathfrak{Q},t)<_{\mathbf{Pth}_{\boldsymbol{\mathcal{A}}^{(1)}}}(\mathfrak{P},s)$, then the following equality holds
\allowdisplaybreaks
\begin{align*}
\mathrm{sc}_{\boldsymbol{\mathcal{B}}^{(1)},\varphi(t)}^{(0,1)}\left(f^{(1)\flat}_{t}\left(\mathfrak{Q}\right)\right) 
&=
f^{(0)\sharp}_{t}\left(\mathrm{sc}_{\boldsymbol{\mathcal{A}}^{(1)},t}^{(0,1)}\left(\mathfrak{Q}\right)\right).
\end{align*}

By Lemma~\ref{LOrdI}, we have that $\mathfrak{P}$ is either~(1) a path of length strictly greater than one containing at least one echelon or~(2) an echelonless path.

If~(1), i.e., if $\mathfrak{P}$ is a path of length strictly greater than one containing at least one echelon, then let $i\in \bb{\mathfrak{P}}$ be the first index for which the one-step subpath $\mathfrak{P}^{i,i}$ of $\mathfrak{P}$ is an echelon. We consider different cases for $i$ according to the cases presented in Definition~\ref{DOrd}.

If $i=0$, we have that the pairs $(\mathfrak{P}^{0,0},s)$ and $(\mathfrak{P}^{1,\bb{\mathfrak{P}}-1},s)$ $\prec_{\mathbf{Pth}_{\boldsymbol{\mathcal{A}}^{(1)}}}$-precede the pair $(\mathfrak{P},s)$. 
The following chain of equalities holds
\allowdisplaybreaks
\begin{align*}
\mathrm{sc}_{\boldsymbol{\mathcal{B}}^{(1)},\varphi(s)}^{(0,1)}\left(
f_{s}^{(1)\flat}\left(
\mathfrak{P}
\right)
\right)&=
\mathrm{sc}_{\boldsymbol{\mathcal{B}}^{(1)},\varphi(s)}^{(0,1)}\left(
f^{(1)\flat}_{s}\left(
\mathfrak{P}^{1,\bb{\mathfrak{P}}-1}
\right)
\circ_{\varphi(s)}^{0\mathbf{Pth}_{\boldsymbol{\mathcal{B}}^{(1)}}}
f^{(1)\flat}_{s}\left(
\mathfrak{P}^{0,0}
\right)
\right)
\tag{1}
\\
&=
\mathrm{sc}_{\boldsymbol{\mathcal{B}}^{(1)},\varphi(s)}^{(0,1)}\left(
f^{(1)\flat}_{s}\left(
\mathfrak{P}^{0,0}
\right)
\right)
\tag{2}
\\
&=
f^{(0)\sharp}_{s}\left(
\mathrm{sc}_{\boldsymbol{\mathcal{A}}^{(1)},s}^{(0,1)}\left(
\mathfrak{P}^{0,0}
\right)
\right)
\tag{3}
\\
&=
f^{(0)\sharp}_{s}\left(
\mathrm{sc}_{\boldsymbol{\mathcal{A}}^{(1)},s}^{(0,1)}\left(
\mathfrak{P}^{1,\bb{\mathfrak{P}}-1}
\circ_{s}^{0\mathbf{Pth}_{\boldsymbol{\mathcal{A}}^{(1)}}}
\mathfrak{P}^{0,0}
\right)
\right)
\tag{4}
\\
&=
f^{(0)\sharp}_{s}\left(
\mathrm{sc}_{\boldsymbol{\mathcal{A}}^{(1)},s}^{(0,1)}\left(
\mathfrak{P}
\right)
\right).
\tag{5}
\end{align*}

The first equality unravels the definition of $f^{(1)\flat}$;
the second equality follows from Proposition~\ref{PPthComp};
the third equality follows by Artinian induction;
the fourth equality follows from Proposition~\ref{PPthComp};
the fifth equality recovers the definition of the path $\mathfrak{P}$.

If $i\neq 0$, we have that the pairs $(\mathfrak{P}^{0,i-1},s)$ and $(\mathfrak{P}^{i,\bb{\mathfrak{P}}-1},s)$ $\prec_{\mathbf{Pth}_{\boldsymbol{\mathcal{A}}^{(1)}}}$-precede the pair $(\mathfrak{P},s)$. 
The following chain of equalities holds
\allowdisplaybreaks
\begin{align*}
\mathrm{sc}_{\boldsymbol{\mathcal{B}}^{(1)},\varphi(s)}^{(0,1)}\left(
f_{s}^{(1)\flat}\left(
\mathfrak{P}
\right)
\right)&=
\mathrm{sc}_{\boldsymbol{\mathcal{B}}^{(1)},\varphi(s)}^{(0,1)}\left(
f^{(1)\flat}_{s}\left(
\mathfrak{P}^{i,\bb{\mathfrak{P}}-1}
\right)
\circ_{\varphi(s)}^{0\mathbf{Pth}_{\boldsymbol{\mathcal{B}}^{(1)}}}
f^{(1)\flat}_{s}\left(
\mathfrak{P}^{0,i-1}
\right)
\right)
\tag{1}
\\
&=
\mathrm{sc}_{\boldsymbol{\mathcal{B}}^{(1)},\varphi(s)}^{(0,1)}\left(
f^{(1)\flat}_{s}\left(
\mathfrak{P}^{0,i-1}
\right)
\right)
\tag{2}
\\
&=
f^{(0)\sharp}_{s}\left(
\mathrm{sc}_{\boldsymbol{\mathcal{A}}^{(1)},s}^{(0,1)}\left(
\mathfrak{P}^{0,1-i}
\right)
\right)
\tag{3}
\\
&=
f^{(0)\sharp}_{s}\left(
\mathrm{sc}_{\boldsymbol{\mathcal{A}}^{(1)},s}^{(0,1)}\left(
\mathfrak{P}^{i,\bb{\mathfrak{P}}-1}
\circ_{s}^{1\mathbf{Pth}_{\boldsymbol{\mathcal{A}}^{(1)}}}
\mathfrak{P}^{0,i-1}
\right)
\right)
\tag{4}
\\
&=
f^{(0)\sharp}_{s}\left(
\mathrm{sc}_{\boldsymbol{\mathcal{A}}^{(1)},s}^{(0,1)}\left(
\mathfrak{P}
\right)
\right).
\tag{5}
\end{align*}

The first equality unravels the definition of $f^{(1)\flat}$;
the second equality follows from Proposition~\ref{PPthComp};
the third equality follows by Artinian induction;
the fourth equality follows from Proposition~\ref{PPthComp};
the fifth equality recovers the definition of the path $\mathfrak{P}$.

If~(2), i.e., if $\mathfrak{P}$ is an echelonless path in $\mathrm{Pth}_{\boldsymbol{\mathcal{A}}^{(1)},s}$, then the conditions for the path extraction algorithm, as stated in Lemma~\ref{LPthExtract}, are fulfilled. Then, by Lemma~\ref{LPthHeadCt}, there exists a unique word $\mathbf{s} \in S^{\star} - \{ \lambda \}$ and a unique operation symbol $\sigma \in \Sigma_{\mathbf{s}, s}$ associated to $\mathfrak{P}$. Let $(\mathfrak{P}_{j})_{j \in \bb{\mathbf{s}}}$ be the family of paths in $\mathrm{Pth}_{\boldsymbol{\mathcal{A}}^{(1)}, s}$ which, in virtue of Lemma~\ref{LPthExtract}, we can extract from $\mathfrak{P}$. Note that, for every $j \in \bb{\mathbf{s}}$, we have that $(\mathfrak{P}_{j}, s_{j}) \prec_{\mathbf{Pth}_{\boldsymbol{\mathcal{A}}^{(1)}}} (\mathfrak{P}, s)$. 

The following chain of equalities holds
\allowdisplaybreaks
\begin{align*}
\mathrm{sc}_{\boldsymbol{\mathcal{B}}^{(1)},\varphi(s)}^{(0,1)}\left(
f_{s}^{(1)\flat}\left(
\mathfrak{P}
\right)
\right)
&=
\mathrm{sc}_{\boldsymbol{\mathcal{B}}^{(1)},\varphi(s)}^{(0,1)}\left(
\sigma^{\mathbf{Pth}_{\boldsymbol{\mathcal{B}}^{(1)}}^{\mathbf{f}^{(1)}(0,1)}}\left(\left(
f_{s_{j}}^{(1)\flat}\left(
\mathfrak{P}_{j}
\right)
\right)_{j\in\bb{\mathbf{s}}}\right)
\right)
\tag{1}
\\
&=
\sigma^{\mathbf{c}_{\mathfrak{d}}^{\ast}(\mathbf{T}_{\Lambda}(Y))}\left(\left(
\mathrm{sc}_{\boldsymbol{\mathcal{B}}^{(1)},\varphi(s_{j})}^{(0,1)}\left(
f_{s_{j}}^{(1)\flat}\left(
\mathfrak{P}_{j}
\right)
\right)
\right)_{j\in\bb{\mathbf{s}}}\right)
\tag{2}
\\
&=
\sigma^{\mathbf{c}_{\mathfrak{d}}^{\ast}(\mathbf{T}_{\Lambda}(Y))}\left(\left(
f_{s_{j}}^{(0)\sharp}\left(
\mathrm{sc}_{\boldsymbol{\mathcal{A}}^{(1)},s_{j}}^{(0,1)}\left(
\mathfrak{P}_{j}
\right)
\right)
\right)_{j\in\bb{\mathbf{s}}}\right)
\tag{3}
\\
&=
f_{s}^{(0)\sharp}\left(
\sigma^{\mathbf{T}_{\Sigma}(X)}\left(\left(
\mathrm{sc}_{\boldsymbol{\mathcal{A}}^{(1)},s_{j}}^{(0,1)}\left(
\mathfrak{P}_{j}
\right)
\right)_{j\in\bb{\mathbf{s}}}\right)
\right)
\tag{4}
\\
&=
f_{s}^{(0)\sharp}\left(
\mathrm{sc}_{\boldsymbol{\mathcal{A}}^{(1)},s}^{(0,1)}\left(
\sigma^{\mathbf{Pth}_{\boldsymbol{\mathcal{A}}^{(1)}}}\left(\left(
\mathfrak{P}_{j}
\right)_{j\in\bb{\mathbf{s}}}\right)
\right)
\right)
\tag{5}
\\
&=
f_{s}^{(0)\sharp}\left(
\mathrm{sc}_{\boldsymbol{\mathcal{A}}^{(1)},s}^{(0,1)}\left(
\mathfrak{P}
\right)
\right).
\tag{6}
\end{align*}

The first equality unravels the definition of $f^{(1)\flat}$;
the second equality follows from the fact that, by Remark~\ref{RDSigmaAlg}, $\mathrm{sc}_{\boldsymbol{\mathcal{B}}^{(1)}, \varphi}^{(0,1)}$ is a $\Sigma$-homomorphism from $\mathbf{Pth}^{\mathbf{f}^{(1)}(0,1)}_{\boldsymbol{\mathcal{B}}^{(1)}}$ to $\mathbf{c}_{\mathfrak{d}}^{\ast}(\mathbf{T}_{\Lambda}(Y))$;
the third equality follows by Artinian induction;
the fourth equality follows from the fact that, following Definition~\ref{DTw0}, $f^{(0)\sharp}$ is a $\Sigma$-homomorphism from $\mathbf{T}_{\Sigma}(X)$ to $\mathbf{c}_{\mathfrak{d}}^{\ast}(\mathbf{T}_{\Lambda}(Y))$;
the fifth equality follows from Proposition~\ref{PPthAlg};
finally, the last equality follows from the definition of $\mathfrak{P}$.

This completes the proof of Item (1).

\textsf{(2)} $\mathrm{tg}^{(0,1)}_{\boldsymbol{\mathcal{B}}^{(1)},\varphi}\circ f^{(1)\flat}=f^{(0)\sharp}\circ \mathrm{tg}^{(0,1)}_{\boldsymbol{\mathcal{A}}^{(1)}}$.

This proof is similar to that of item (1).

\textsf{(3)} $f^{(1)\flat} \circ \mathrm{ip}_{\boldsymbol{\mathcal{A}}^{(1)}}^{(1, 0)\sharp} = \mathrm{ip}_{\boldsymbol{\mathcal{B}}^{(1)}, \varphi}^{(1, 0)\sharp} \circ f^{(0)\sharp}$.

This property follows directly from the definition of the mapping $f^{(1)\flat}$.

\textsf{(4)} $f^{(1)\flat} \circ \mathrm{ech}^{(1,\mathcal{A}^{(1)})}_{\boldsymbol{\mathcal{A}}^{(1)}} = f^{(1)}$. 

This property follows directly from the definition of the mapping $f^{(1)\flat}$.

This completes the proof.
\end{proof}

We next prove that, for a first-order morphism, its path extension mapping introduced in Proposition~\ref{PPthExt} is a $\Sigma$-homomorphism.

\begin{proposition}
\label{PPthExtHom}
Let $\mathbf{f}^{(1)}=(\varphi, c, (f^{(i)})_{i\in 2})$ be a first-order morphism from $\boldsymbol{\mathcal{A}}^{(1)}$ to $\boldsymbol{\mathcal{B}}^{(1)}$. Then the \emph{path extension mapping of $f^{(1)}$}, $f^{(1)\flat}\colon \mathrm{Pth}_{\boldsymbol{\mathcal{A}^{(1)}}}\mor \mathrm{Pth}_{\boldsymbol{\mathcal{B}^{(1)},\varphi}}$ introduced in Proposition~\ref{PPthExt} is a $\Sigma$-homomorphism
$$
f^{(1)\flat}\colon
\mathrm{Pth}^{(0,1)}_{\boldsymbol{\mathcal{A}}^{(1)}}
\mor 
\mathrm{Pth}^{\mathbf{f}^{(1)}(0,1)}_{\boldsymbol{\mathcal{B}}^{(1)}}.
$$
\end{proposition}

\begin{proof}
Let $(\mathbf{s},s)$ be an element of $S^{\star}\times S$ and $\sigma$ an operation symbol in $\Sigma_{\mathbf{s},s}$.

If $\mathbf{s}=\lambda$, then the following chain of equalities holds
\allowdisplaybreaks
\begin{align*}
f_{s}^{(1)\flat}\left(
\sigma^{\mathbf{Pth}_{\boldsymbol{\mathcal{A}}^{(1)}}^{(0,1)}}
\right)
&=
f_{s}^{(1)\flat}\left(
\mathrm{ip}^{(1,0)\sharp}_{\boldsymbol{\mathcal{A}}^{(1)},s}\left(
\sigma^{\mathbf{T}_{\Sigma}(X)}
\right)
\right)
\tag{1}
\\
&=
\mathrm{ip}^{(1,0)\sharp}_{\boldsymbol{\mathcal{B}}^{(1)},\varphi(s)}\left(
f_{s}^{(0)\sharp}\left(
\sigma^{\mathbf{T}_{\Sigma}(X)}
\right)
\right)
\tag{2}
\\
&=
\mathrm{ip}_{\boldsymbol{\mathcal{B}}^{(1)},\varphi(s)}^{(1,0)\sharp}\left(
\sigma^{\mathbf{c}_{\mathfrak{d}}^{\ast}(\mathbf{T}_{\Lambda}(Y))}
\right)
\tag{3}
\\
&=
\sigma^{\mathbf{Pth}_{\boldsymbol{\mathcal{B}}^{(1)}}^{\mathbf{f}^{(1)}(0,1)}}.
\tag{4}
\end{align*}

The first equality follows from the fact that, for a constant $\sigma$ in $\Sigma_{\lambda, s}$, then, following Remark~\ref{RConsSigma}, the interpretation of $\sigma$ in $\mathbf{Pth}_{\boldsymbol{\mathcal{A}}^{(1)}}$ is given by $\sigma^{\mathbf{Pth}_{\boldsymbol{\mathcal{A}}^{(1)}}} = \mathrm{ip}_{s}^{(1,0)\sharp}(\sigma^{\mathbf{T}_{\Sigma}(X)})$;
the second equality unravels the definition of $f^{(1)\flat}$ at a $(1,0)$-identity path;
The third equality follows from the fact that, according to Proposition~\ref{PPropUniv}, $f^{(0)\sharp}$ is a $\Sigma$-homomorphism from $\mathbf{T}_{\Sigma}(X)$ to $\mathbf{c}_{\mathfrak{d}}^{\ast}(\mathbf{T}_{\Lambda}(Y))$;
finally, the last equality from the fact that, according to Remark~\ref{RDSigmaAlg}, $\mathrm{ip}_{\boldsymbol{\mathcal{B}}^{(1)},\varphi(s)}^{(1,0)\sharp}$ is a $\Sigma$-homomorphism from $\mathbf{c}_{\mathfrak{d}}^{\ast}(\mathbf{T}_{\Lambda}(Y))$ to $\mathbf{Pth}_{\boldsymbol{\mathcal{B}}^{(1)}}^{\mathbf{f}^{(1)}(0,1)}$.

We now consider the case in which $\mathbf{s}\neq\lambda$. Let $(\mathfrak{P}_{j})_{j\in\bb{\mathbf{s}}}$ be a family of paths in $\mathrm{Pth}_{\boldsymbol{\mathcal{A}}^{(1)},\mathbf{s}}$. We consider different cases according to the nature of the family $(\mathfrak{P}_{j})_{j\in\bb{\mathbf{s}}}$. It could be the case that either~(1), for every $j\in\bb{\mathbf{s}}$, $\mathfrak{P}_{j}$ is a $(1,0)$-identity path or~(2), there exists an index $j\in\bb{\mathbf{s}}$ for which $\mathfrak{P}_{j}$ is not an $(1,0)$-identity path.

If~(1), then for every $j\in\bb{\mathbf{s}}$, $\mathfrak{P}_{j}$ is equal to $\mathrm{ip}_{\boldsymbol{\mathcal{A}}^{(1)},s_{j}}^{(1,0)\sharp}(P_{j})$ for some term $P_{j}$ in $\T_{\Sigma}(X)_{s_{j}}$. In this case, the following chain of equalities holds
\allowdisplaybreaks
\begin{align*}
f_{s}^{(1)\flat}\left(
\sigma^{\mathbf{Pth}_{\boldsymbol{\mathcal{A}}^{(1)}}^{(0,1)}}\left(\left(
\mathfrak{P}_{j}
\right)_{j\in\bb{\mathbf{s}}}\right)
\right)&=
f_{s}^{(1)\flat}\left(
\sigma^{\mathbf{Pth}_{\boldsymbol{\mathcal{A}}^{(1)}}^{(0,1)}}\left(\left(
\mathrm{ip}_{\boldsymbol{\mathcal{A}}^{(1)},s_{j}}^{(1,0)\sharp}\left(
P_{j}
\right)
\right)_{j\in\bb{\mathbf{s}}}\right)
\right)
\tag{1}
\\
&=
f_{s}^{(1)\flat}\left(
\mathrm{ip}_{\boldsymbol{\mathcal{A}}^{(1)},s}^{(1,0)\sharp}\left(
\sigma^{\mathbf{T}_{\Sigma}(X)}\left(\left(
P_{j}
\right)_{j\in\bb{\mathbf{s}}}\right)
\right)
\right)
\tag{2}
\\
&=
\mathrm{ip}_{\boldsymbol{\mathcal{B}}^{(1)},\varphi(s)}^{(1,0)\sharp}\left(
f_{s}^{(0)\sharp}\left(
\sigma^{\mathbf{T}_{\Sigma}(X)}\left(\left(
P_{j}
\right)_{j\in\bb{\mathbf{s}}}\right)
\right)
\right)
\tag{3}
\\
&=
\mathrm{ip}_{\boldsymbol{\mathcal{B}}^{(1)},\varphi(s)}^{(1,0)\sharp}\left(
\sigma^{\mathbf{c}_{\mathfrak{d}}^{\ast}(\mathbf{T}_{\Lambda}(Y))}\left(\left(
f_{s_{j}}^{(0)\sharp}\left(
P_{j}
\right)\right)_{j\in\bb{\mathbf{s}}}
\right)
\right)
\tag{4}
\\
&=
\sigma^{\mathbf{Pth}_{\boldsymbol{\mathcal{B}}^{(1)}}^{\mathbf{f}^{(1)}(0,1)}}\left(\left(
\mathrm{ip}_{\boldsymbol{\mathcal{B}}^{(1)},\varphi(s_{j})}^{(1,0)\sharp}\left(
f_{s_{j}}^{(0)\sharp}\left(
P_{j}
\right)\right)
\right)_{j\in\bb{\mathbf{s}}}\right)
\tag{5}
\\
&=
\sigma^{\mathbf{Pth}_{\boldsymbol{\mathcal{B}}^{(1)}}^{\mathbf{f}^{(1)}(0,1)}}\left(\left(
{f}_{s_{j}}^{(1)\flat}\left(
\mathrm{ip}_{\boldsymbol{\mathcal{A}}^{(1)},s_{j}}^{(1,0)\sharp}\left(
P_{j}
\right)
\right)
\right)_{j\in\bb{\mathbf{s}}}\right)
\tag{6}
\\
&=
\sigma^{\mathbf{Pth}_{\boldsymbol{\mathcal{B}}^{(1)}}^{\mathbf{f}^{(1)}(0,1)}}\left(\left(
{f}_{s_{j}}^{(1)\flat}\left(
\mathfrak{P}_{j}
\right)
\right)_{j\in\bb{\mathbf{s}}}\right).
\tag{7}
\end{align*}

The first equality follows from the fact that, for every $j\in\bb{\mathbf{s}}$, $\mathfrak{P}_{j}$ is a $(1,0)$-identity path;
the first equality follows from the fact that, by Proposition~\ref{PIpHom}, $\mathrm{ip}_{\boldsymbol{\mathcal{A}}^{(1)}}^{(1,0)\sharp}$ is a $\Sigma$-homomorphism from $\mathbf{T}_{\Sigma}(X)$ to $\mathbf{Pth}_{\boldsymbol{\mathcal{A}}^{(1)}}^{(0,1)}$;
the third equality unravels the definition of $f^{(1)\flat}$ at a $(1,0)$-identity path;
the fourth equality follows from the fact that, by Proposition~\ref{PPropUniv}, $f^{(0)\sharp}$ is a $\Sigma$-homomorphism from $\mathbf{T}_{\Sigma}(X)$ to $\mathbf{c}_{\mathfrak{d}}^{\ast}(\mathbf{T}_{\Lambda}(Y)$;
the fifth equality follows from the fact that, by Remark~\ref{RDSigmaAlg}, $\mathrm{ip}^{(1,0)\sharp}_{\boldsymbol{\mathcal{B}}^{(1)}, \varphi}$ is a $\Sigma$-homomorphism from $\mathbf{c}_{\mathfrak{d}}^{\ast}(\mathbf{T}^{\Lambda}(Y))$ to $\mathbf{Pth}_{\boldsymbol{\mathcal{B}}^{(1)}}^{\mathbf{f}^{(1)}(0,1)}$;
the sixth equality recovers the definition of $f^{(1)\flat}$ at a $(1,0)$-identity path;
finally, the last equality recovers, for every $j\in\bb{\mathbf{s}}$, the definition of $\mathfrak{P}_{j}$ as a $(1,0)$-identity path.

This proves Case~(1).

If~(2), i.e., if there exists some index $j\in\bb{\mathbf{s}}$ for which $\mathfrak{P}_{j}$ is not an $(1,0)$-identity path then, according to Corollary~\ref{CPthWB}, $\sigma^{\mathbf{Pth}_{\boldsymbol{\mathcal{A}}^{(1)}}^{(0,1)}}((\mathfrak{P}_{j})_{j\in\bb{\mathbf{s}}})$ is an echelonless path. Moreover, according to Proposition~\ref{PRecov}, the path extraction algorithm from Lemma~\ref{LPthExtract}, applied to it retrieves the original family $(\mathfrak{P}_{j})_{j\in\bb{\mathbf{s}}}$. Then, by definition of the $S$-sorted mapping $f^{(1)\flat}$,
$$
f^{(1)\flat}_{s}\left(
\sigma^{\mathbf{Pth}_{\boldsymbol{\mathcal{A}}^{(1)}}^{(0,1)}}\left(\left(
\mathfrak{P}_{j}
\right)_{j\in\bb{\mathbf{s}}}\right)
\right)
=
\sigma^{\mathbf{Pth}_{\boldsymbol{\mathcal{B}}^{(1)}}^{\mathbf{f}^{(1)}(0,1)}}\left(\left(
f^{(1)\flat}_{s_{j}}\left(
\mathfrak{P}_{j}
\right)
\right)_{j\in\bb{\mathbf{s}}}\right)
$$

Case~(2) follows. 

This finishes the proof.
\end{proof}			
\chapter{Derived partial $\Sigma^{\boldsymbol{\mathcal{A}}^{(1)}}$-algebras}\label{S3B}

In this chapter we define a structure of partial $\Sigma^{\boldsymbol{\mathcal{A}}^{(1)}}$-algebra on $\mathrm{Pth}_{\boldsymbol{\mathcal{B}}^{(1)}, \varphi}$ and $[ \mathrm{Pth}_{\boldsymbol{\mathcal{B}}^{(1)}} ]_{\varphi}$, which we will denote by $\mathbf{Pth}_{\boldsymbol{\mathcal{B}}^{(1)}}^{\mathbf{f}^{(1)}}$ and $[ \mathbf{Pth}_{\boldsymbol{\mathcal{B}}^{(1)}}^{\mathbf{f}^{(1)}} ]$, respectively.
We show that the path extension mapping is not a $\Sigma^{\boldsymbol{\mathcal{A}}^{(1)}}$-homomorphism.
However, we show that the projection mapping $\mathrm{pr}^{\mathrm{Ker}(\mathrm{CH}^{(1)})}_{\boldsymbol{\mathcal{B}}^{(1)}, \varphi}$ is a $\Sigma^{\boldsymbol{\mathcal{A}}^{(1)}}$-homomorphism.
We next define a structure of partial $\Sigma^{\boldsymbol{\mathcal{A}}^{(1)}}$-algebras on $\mathrm{PT}_{\boldsymbol{\mathcal{B}}^{(1)}, \varphi}$ and $[\mathrm{PT}_{\boldsymbol{\mathcal{B}}^{(1)}}]_{\varphi}$, which we will denote by $\mathbf{PT}_{\boldsymbol{\mathcal{B}}^{(1)}}^{\mathbf{f}^{(1)}}$ and $[ \mathbf{PT}_{\boldsymbol{\mathcal{B}}^{(1)}}^{\mathbf{f}^{(1)}} ]$, respectively.
We then show that the projection mapping $\mathrm{pr}^{\Theta^{[1]}}_{\boldsymbol{\mathcal{B}}^{(1)}, \varphi}$ is a $\Sigma^{\boldsymbol{\mathcal{A}}^{(1)}}$-homomorphism.
Moreover, we show that $[ \mathbf{Pth}_{\boldsymbol{\mathcal{B}}^{(1)}}^{\mathbf{f}^{(1)}} ]$ and $[ \mathbf{PT}_{\boldsymbol{\mathcal{B}}^{(1)}}^{\mathbf{f}^{(1)}} ]$ are isomorphic partial $\Sigma^{\boldsymbol{\mathcal{A}}^{(1)}}$-algebras. 
Finally, we define a structure of partial $\Sigma^{\boldsymbol{\mathcal{A}}^{(1)}}$-algebra on $\T_{\boldsymbol{\mathcal{E}}^{\boldsymbol{\mathcal{B}}^{(1)}}}(\mathbf{Pth}_{\boldsymbol{\mathcal{B}}^{(1)}})_{\varphi}$, which we will denote by $\mathbf{T}_{\boldsymbol{\mathcal{E}}^{\boldsymbol{\mathcal{B}}^{(1)}}}^{\mathbf{f}^{(1)}}(\mathbf{Pth}_{\boldsymbol{\mathcal{B}}^{(1)}})$ and we prove that the partial $\Sigma^{\boldsymbol{\mathcal{A}}^{(1)}}$-algebras $[ \mathbf{Pth}_{\boldsymbol{\mathcal{B}}^{(1)}}^{\mathbf{f}^{(1)}} ]$ and $\mathbf{T}_{\boldsymbol{\mathcal{E}}^{\boldsymbol{\mathcal{B}}^{(1)}}}^{\mathbf{f}^{(1)}}(\mathbf{Pth}_{\boldsymbol{\mathcal{B}}^{(1)}})$ are also isomorphic.

\section{A structure of partial $\Sigma^{\boldsymbol{\mathcal{A}}^{(1)}}$-algebra on $\mathrm{Pth}_{\boldsymbol{\mathcal{B}}^{(1)}, \varphi}$}

We next show that, given a first-order morphism $\mathbf{f}^{(1)}$, $\mathrm{Pth}_{\boldsymbol{\mathcal{B}}^{(1)}, \varphi}$ is equipped in a natural way with a structure of partial $\Sigma^{\boldsymbol{\mathcal{A}}^{(1)}}$-algebra that we denote by $\mathbf{Pth}_{\boldsymbol{\mathcal{B}}^{(1)}}^{\mathbf{f}^{(1)}}$.

\begin{proposition}
\label{PPthBCatAlg}
Let $\mathbf{f}^{(1)}=(\varphi, c, (f^{(i)})_{i\in 2})$ be a first-order morphism from $\boldsymbol{\mathcal{A}}^{(1)}$ to $\boldsymbol{\mathcal{B}}^{(1)}$. Then the $S$-sorted set $\mathrm{Pth}_{\boldsymbol{\mathcal{B}}^{(1)}, \varphi}$ is equipped, in a natural way, with a structure of partial $\Sigma^{\boldsymbol{\mathcal{A}}^{(1)}}$-algebra.
\end{proposition}

\begin{proof}
Let us denote by $\mathbf{Pth}_{\boldsymbol{\mathcal{B}}^{(1)}}^{\mathbf{f}^{(1)}}$ the $\Sigma^{\boldsymbol{\mathcal{A}}^{(1)}}$-algebra defined as follows:

\textsf{(1)}
The underlying $S$-sorted set of $\mathbf{Pth}_{\boldsymbol{\mathcal{B}}^{(1)}}^{\mathbf{f}^{(1)}}$ is $\mathrm{Pth}_{\boldsymbol{\mathcal{B}}^{(1)}, \varphi}=(\mathrm{Pth}_{\boldsymbol{\mathcal{B}}^{(1)}, \varphi(s)})_{s\in S}$.

\textsf{(2)}
For every $(\mathbf{s}, s)\in S^{\ast}\times S$ and every operation symbol $\sigma\in\Sigma_{\mathbf{s},s}$, the operation $\sigma^{\mathbf{Pth}_{\boldsymbol{\mathcal{B}}^{(1)}}^{\mathbf{f}^{(1)}}}$ is given by the interpretation of $\sigma$ in the $\Sigma$-algebra $\mathbf{Pth}_{\boldsymbol{\mathcal{B}}^{(1)}}^{\mathbf{f}^{(1)}(0,1)}$ introduced in Remark~\ref{RDSigmaAlg}.

\textsf{(3)}
For every $s\in S$ and every $\mathfrak{p}\in\mathcal{A}_{s}^{(1)}$, the constant $\mathfrak{p}^{\mathbf{Pth}_{\boldsymbol{\mathcal{B}}^{(1)}}^{\mathbf{f}^{(1)}}}$ is given by the image of $\mathfrak{p}$ under the mapping $f^{(1)}$, i.e., we set
$$
\mathfrak{p}^{\mathbf{Pth}_{\boldsymbol{\mathcal{B}}^{(1)}}^{\mathbf{f}^{(1)}}}
=
f^{(1)}_{s}(\mathfrak{p})
\in
\mathrm{Pth}_{\boldsymbol{\mathcal{B}}^{(1)}, \varphi(s)}.
$$

\textsf{(4)}
For every $s\in S$, the interpretations of the operations $\mathrm{sc}_{s}^{0}$ and $\mathrm{tg}_{s}^{0}$ are given by
\begin{align*}
\mathrm{sc}^{0\mathbf{Pth}_{\boldsymbol{\mathcal{B}}^{(1)}}^{\mathbf{f}^{(1)}}}_{s}
&=
\mathrm{sc}^{0\mathbf{Pth}_{\boldsymbol{\mathcal{B}}}^{(1)}}_{\varphi(s)}
&\mbox{and}&&
\mathrm{tg}^{0\mathbf{Pth}_{\boldsymbol{\mathcal{B}}^{(1)}}^{\mathbf{f}^{(1)}}}_{s}
&=
\mathrm{tg}^{0\mathbf{Pth}_{\boldsymbol{\mathcal{B}}}^{(1)}}_{\varphi(s)}
\end{align*}
That is, teir interpretations are given by the interpretations of $\mathrm{sc}^{0}_{\varphi(s)}$ and $\mathrm{tg}^{0}_{\varphi(s)}$ in $\mathbf{Pth}_{\boldsymbol{\mathcal{B}}^{(1)}}$ that, we recall, was introduced in Proposition~\ref{PPthCatAlg}.

\textsf{(5)}
Similarly, for every $s \in S$,  the interpretation of the partial binary operation $\circ_{s}^{0}$ is given by
$$
\circ_{s}^{0\mathbf{Pth}_{\boldsymbol{\mathcal{B}}^{(1)}}^{\mathbf{f}^{(1)}}}
=
\circ_{\varphi(s)}^{0\mathbf{Pth}_{\boldsymbol{\mathcal{B}}^{(1)}}}.
$$
That is, its interpretation is given by the interpretation of $\circ_{\varphi(s)}^{0}$ in $\mathbf{Pth}_{\boldsymbol{\mathcal{B}}^{(1)}}$ that, we recall, was introduced in Proposition~\ref{PPthCatAlg}.

This completes the definition of the partial $\Sigma^{\boldsymbol{\mathcal{A}}^{(1)}}$-algebra $\mathbf{Pth}_{\boldsymbol{\mathcal{B}}^{(1)}}^{\mathbf{f}^{(1)}}$.
\end{proof}

Despite the relations of the path extension mapping with the source, target, and identity path mappings, it cannot be inferred that it is a $\Sigma^{\boldsymbol{\mathcal{A}}^{(1)}}$-homomorphism. This is largely due to the relationship between the $0$-composition of paths and the Artinian order used to define the mapping $f^{(1)\flat}$. Below we present a counterexample to the assertion that $f^{(1)\flat}$ is a $\Sigma^{\boldsymbol{\mathcal{A}}^{(1)}}$-homomorphism.

\begin{remark}
Let $\mathbf{f}^{(1)} = (c, \varphi, (f^{(i)})_{i \in 2})$ be a first-order morphism from $\boldsymbol{\mathcal{A}}^{(1)}$ to $\boldsymbol{\mathcal{B}}^{(1)}$. The path extnesion mapping $f^{(1)\flat}$ from $\mathrm{Pth}_{\boldsymbol{\mathcal{A}}^{(1)}}$ to $\mathrm{Pth}_{\boldsymbol{\mathcal{B}}^{(1)}, \varphi}$ is not necessarily a $\Sigma^{\boldsymbol{\mathcal{A}}^{(1)}}$-homomorphism from $\mathbf{Pth}_{\boldsymbol{\mathcal{A}}^{(1)}}$ to $\mathbf{Pth}_{\boldsymbol{\mathcal{B}}^{(1)}}^{\mathbf{f}^{(1)}}$.

Let us consider a path of the form $\mathfrak{Q}\circ_{s}^{0\mathbf{Pth}_{\boldsymbol{\mathcal{A}}^{(1)}}}\mathfrak{P}$ and let us suppose that it is echelonless. Then, regarding the paths $\mathfrak{Q}$ and $\mathfrak{P}$ we have that

\begin{itemize}
    \item[(i)] $\mathfrak{P}$ is an echelonless path;
    \item[(ii)] $\mathfrak{Q}$ is an echelonless path.
\end{itemize}

In this case, by Lemma~\ref{LPthHeadCt} there exists a unique word $\mathbf{s}\in S^{\star}-\{\lambda\}$ and a unique operation symbol $\sigma\in\Sigma_{\mathbf{s},s}$ such that the family of translations ocurring in $\mathfrak{Q}\circ_{s}^{0\mathbf{Pth}_{\boldsymbol{\mathcal{A}}^{(1)}}}\mathfrak{P}$ is a family of translations of type $\sigma$.

Since $\mathfrak{Q}\circ_{s}^{0\mathbf{Pth}_{\boldsymbol{\mathcal{A}}^{(1)}}}\mathfrak{P}$ is associated to the operation symbol $\sigma$, in virtue of Lemma~\ref{LPthHeadCt} all translations appearing in $\mathfrak{Q}$ and $\mathfrak{P}$ have the same type, i.e., the paths $\mathfrak{Q}$ and $\mathfrak{P}$ are also associated to this same operation symbol $\sigma$.

Let $((\mathfrak{Q}\circ_{s}^{0\mathbf{Pth}_{\boldsymbol{\mathcal{A}}^{(1)}}}\mathfrak{P}) _{j})_{j\in\bb{\mathbf{s}}})$ be the family of paths we can extract from $\mathfrak{Q}\circ_{s}^{0\mathbf{Pth}_{\boldsymbol{\mathcal{A}}^{(1)}}}\mathfrak{P}$ in virtue of Lemma\ref{LPthExtract}. Then, the value of the mapping $f^{(1)\flat}$ at $\mathfrak{Q}\circ_{s}^{0\mathbf{Pth}_{\boldsymbol{\mathcal{A}}^{(1)}}}\mathfrak{P}$ is given by

$$
f_{s}^{(1)\flat}\left(
\mathfrak{Q}\circ_{s}^{0\mathbf{Pth}_{\boldsymbol{\mathcal{A}}^{(1)}}}\mathfrak{P}
\right)
=
\sigma^{\mathbf{Pth}_{\boldsymbol{\mathcal{B}}}^{\mathbf{f}^{(1)}(0,1)}(0,1)}\left(
\left(
f_{s_{j}}^{(1)\flat}\left(
(\mathfrak{Q}\circ_{s}^{0\mathbf{Pth}_{\boldsymbol{\mathcal{A}}^{(1)}}}\mathfrak{P})_{j}
\right)
\right)_{j\in\bb{\mathbf{s}}}
\right).
$$

Moreover, note that
\begin{itemize}
    \item[(i)] $\mathfrak{P}$ is an echelonless path associated to the operation symbol $\sigma$.
\end{itemize}

Let $(\mathfrak{P}_{j})_{j\in\bb{\mathbf{s}}}$ be the family of paths we can extract from $\mathfrak{P}$ in virtue of Lemma~\ref{LPthExtract}. Then, according to Proposition~\ref{PPthExt} the value of the mapping $f^{(1)\flat}$ at $\mathfrak{P}$ is given by
$$
f_{s}^{(1)\flat}(\mathfrak{P})
=
\sigma^{\mathbf{Pth}_{\boldsymbol{\mathcal{B}}^{(1)}}^{\mathbf{f}^{(1)}(0,1)}}
\left(
\left(
f_{s_{j}}^{(1)\flat}(\mathfrak{P}_{j})
\right)_{j\in\bb{\mathbf{s}}}
\right).
$$

On the other hand

\begin{itemize}
    \item[(ii)] $\mathfrak{Q}$ is an echelonless path associated to the operation symbol $\sigma$.
\end{itemize}

Let $(\mathfrak{Q}_{j})_{j\in\bb{\mathbf{s}}}$ be the family of paths we can extract from $\mathfrak{Q}$ in virtue of Lemma~\ref{LPthExtract}. Then, according to Proposition~\ref{PPthExt} the value of the mapping $f^{(1)\flat}$ at $\mathfrak{Q}$ is given by
$$
f_{s}^{(1)\flat}(\mathfrak{Q})
=
\sigma^{\mathbf{Pth}_{\boldsymbol{\sigma}^{(1)}}^{\mathbf{f}^{(1)}(0,1)}}
\left(
\left(
f_{s_{j}}^{(1)\flat}(\mathfrak{Q}_{j})
\right)_{j\in\bb{\mathbf{s}}}
\right).
$$

All in all, in order to $f^{(1)\flat}$ to be a $\Sigma^{\boldsymbol{\sigma}^{(1)}}$-homomorphism, the equality
\begin{multline*}
\sigma^{\mathbf{Pth}_{\boldsymbol{\mathcal{B}}^{(1)}}^{\mathbf{f}^{(1)}(0,1)}}\left(\left(f_{s_{j}}^{(1)\flat}\left((\mathfrak{Q}\circ_{s}^{0\mathbf{Pth}_{\boldsymbol{\mathcal{A}}^{(1)}}}\mathfrak{P})_{j}\right)\right)_{j\in\bb{\mathbf{s}}}\right)=
\\
\sigma^{\mathbf{Pth}_{\boldsymbol{\mathcal{B}}^{(1)}}^{\mathbf{f}^{(1)}(0,1)}}\left(\left(f_{s_{j}}^{(1)\flat}(\mathfrak{P}_{j})\right)_{j\in\bb{\mathbf{s}}}\right)\circ_{\varphi(s)}^{0\mathbf{Pth}_{\boldsymbol{\mathcal{B}}^{(1)}}}\sigma^{\mathbf{Pth}_{\boldsymbol{\mathcal{B}}^{(1)}}^{\mathbf{f}^{(1)}(0,1)}}\left(\left(f_{s_{j}}^{(1)\flat}(\mathfrak{Q}_{j})\right)_{j\in\bb{\mathbf{s}}}\right).
\end{multline*}
should hold. However, consider the many-sorted rewriting system $((S,\Sigma,X),\boldsymbol{\mathcal{A}}^{(1)})$ and the morphism $((\mathrm{id}^{S},\mathrm{id}^{\Sigma},\eta^{X}),\mathrm{ech}^{(1,\mathcal{A}^{(1)})})$ from $\boldsymbol{\mathcal{A}}^{(1)}$ to $\boldsymbol{\mathcal{A}}^{(1)}$. Thus, the above equation is simplified to
\begin{multline*}
\sigma^{\mathbf{Pth}_{\boldsymbol{\mathcal{A}}^{(1)}}^{(0,1)}}\left(\left((\mathfrak{Q}\circ_{s}^{0\mathbf{Pth}_{\boldsymbol{\mathcal{A}}^{(1)}}}\mathfrak{P})_{j}\right)_{j\in\bb{\mathbf{s}}}\right)=
\\
\sigma^{\mathbf{Pth}_{\boldsymbol{\mathcal{A}}^{(1)}}^{(0,1)}}\left((\mathfrak{P}_{j})_{j\in\bb{\mathbf{s}}}\right)\circ_{s}^{0\mathbf{Pth}_{\boldsymbol{\mathcal{A}}^{(1)}}}\sigma^{\mathbf{Pth}_{\boldsymbol{\mathcal{A}}^{(1)}}^{(0,1)}}\left((\mathfrak{Q}_{j})_{j\in\bb{\mathbf{s}}}\right).
\end{multline*}

which is not necessarily true.
\end{remark}

Now, given a first-order morphism $\mathbf{f}^{(1)}$, $[\mathrm{Pth}_{\boldsymbol{\mathcal{B}}^{(1)}}]_{\varphi}$ is equipped in a natural way with a structure of partial $\Sigma^{\boldsymbol{\mathcal{A}}^{(1)}}$-algebra that we denote by $[\mathbf{Pth}_{\boldsymbol{\mathcal{B}}^{(1)}}^{\mathbf{f}^{(1)}}]$. To end this section we show that $\mathrm{pr}^{\mathrm{Ker}(\mathrm{CH}^{(1)})}_{\boldsymbol{\mathcal{B}}^{(1)}, \varphi}$ is a surjective $\Sigma^{\boldsymbol{\mathcal{A}}^{(1)}}$-homomorphism from $\mathbf{Pth}_{\boldsymbol{\mathcal{B}}^{(1)}}^{\mathbf{f}^{(1)}}$ to $[\mathbf{Pth}_{\boldsymbol{\mathcal{B}}^{(1)}}^{\mathbf{f}^{(1)}}]$.

\begin{proposition}
\label{PQPthBCatAlg}
Let $\mathbf{f}^{(1)}=(\varphi, c, (f^{(i)})_{i\in 2})$ be a first-order morphism from $\boldsymbol{\mathcal{A}}^{(1)}$ to $\boldsymbol{\mathcal{B}}^{(1)}$. Then the $S$-sorted set $[\mathrm{Pth}_{\boldsymbol{\mathcal{B}}^{(1)}}]_{\varphi}$ is equipped, in a natural way, with a structure of partial $\Sigma^{\boldsymbol{\mathcal{A}}^{(1)}}$-algebra.
\end{proposition}

\begin{proof}
Let us denote by $[\mathbf{Pth}_{\boldsymbol{\mathcal{B}}^{(1)}}^{\mathbf{f}^{(1)}}]$ the $\Sigma^{\boldsymbol{\mathcal{A}}^{(1)}}$-algebra defined as follows:

\textsf{(1)}
The underlying $S$-sorted set of $[\mathbf{Pth}_{\boldsymbol{\mathcal{B}}^{(1)}}^{\mathbf{f}^{(1)}}]$ is $[\mathrm{Pth}_{\boldsymbol{\mathcal{B}}^{(1)}}]_{\varphi}=([\mathrm{Pth}_{\boldsymbol{\mathcal{B}}^{(1)}}]_{\varphi(s)})_{s\in S}$.

\textsf{(2)}
For every $(\mathbf{s}, s)\in S^{\ast}\times S$ and every operation symbol $\sigma\in\Sigma_{\mathbf{s},s}$, the operation $\sigma^{[\mathbf{Pth}_{\boldsymbol{\mathcal{B}}^{(1)}}^{\mathbf{f}^{(1)}}]}$ is given by the interpretation of $\sigma$ in the $\Sigma$-algebra $[\mathbf{Pth}_{\boldsymbol{\mathcal{B}}^{(1)}}^{\mathbf{f}^{(1)}(0,1)}]$ introduced in Rermark~\ref{RDSigmaAlg}.

\textsf{(3)}
For every $s\in S$ and every $\mathfrak{p}\in\mathcal{A}_{s}^{(1)}$, the constant $\mathfrak{p}^{[\mathbf{Pth}_{\boldsymbol{\mathcal{B}}^{(1)}}^{\mathbf{f}^{(1)}}]}$ is given by the path class of the image of $\mathfrak{p}$ under the mapping $f^{(1)}$, i.e., we set
$$
\mathfrak{p}^{[\mathbf{Pth}_{\boldsymbol{\mathcal{B}}^{(1)}}^{\mathbf{f}^{(1)}}]}
=
\left[
f^{(1)}_{s}\left(
\mathfrak{p}
\right)
\right]_{\varphi(s)}
\in
[
\mathrm{Pth}_{\boldsymbol{\mathcal{B}}^{(1)}}
]_{\varphi(s)}
$$.

\textsf{(4)}
For every $s\in S$, the interpretations of the operations $\mathrm{sc}_{s}^{0}$ and $\mathrm{tg}_{s}^{0}$ are given by
\allowdisplaybreaks
\begin{align*}
\mathrm{sc}_{s}^{0[\mathbf{Pth}_{\boldsymbol{\mathcal{B}}^{(1)}}^{\mathbf{f}^{(1)}}]}
&=
\mathrm{sc}_{\varphi(s)}^{0[\mathbf{Pth}_{\boldsymbol{\mathcal{B}}^{(1)}}]}
&&\mbox{and}&
\mathrm{tg}_{s}^{0[\mathbf{Pth}_{\boldsymbol{\mathcal{B}}^{(1)}}^{\mathbf{f}^{(1)}}]}
&=
\mathrm{tg}_{\varphi(s)}^{0[\mathbf{Pth}_{\boldsymbol{\mathcal{B}}^{(1)}}]}.
\end{align*}
That is, their interpretations are given by the interpretations of $\mathrm{sc}_{\varphi(s)}^{0}$ and $\mathrm{tg}_{\varphi(s)}^{0}$ in $[\mathbf{Pth}_{\boldsymbol{\mathcal{B}}^{(1)}}]$ that, we recall, was introduced in Proposition~\ref{PCHCatAlg}.

\textsf{(5)}
Similarly, for every $s \in S$,  the interpretation of the partial binary operation $\circ_{s}^{0}$ is given by
$$
\circ_{s}^{0[\mathbf{Pth}_{\boldsymbol{\mathcal{B}}^{(1)}}^{\mathbf{f}^{(1)}}]}
=
\circ_{\varphi(s)}^{0[\mathbf{Pth}_{\boldsymbol{\mathcal{B}}^{(1)}}]}.
$$
That is, its interpretation is given by the interpretation of $\circ_{\varphi(s)}^{0}$ in $[\mathbf{Pth}_{\boldsymbol{\mathcal{B}}^{(1)}}]$ that, we recall, was introduced in Proposition~\ref{PCHCatAlg}.

This completes the definition of the partial $\Sigma^{\boldsymbol{\mathcal{A}}^{(1)}}$-algebra $[\mathbf{Pth}_{\boldsymbol{\mathcal{B}}^{(1)}}^{\mathbf{f}^{(1)}}]$.
\end{proof}

\begin{proposition}
\label{PKerCHBCatHom}
Let $\mathbf{f}^{(1)}=(\varphi, c, (f^{(i)})_{i\in 2})$ be a first-order morphism from $\boldsymbol{\mathcal{A}}^{(1)}$ to $\boldsymbol{\mathcal{B}}^{(1)}$. Then the mapping 
$$
\mathrm{pr}_{\boldsymbol{\mathcal{B}}^{(1)}, \varphi}^{\mathrm{Ker}(\mathrm{CH}^{(1)})}
\colon 
\mathrm{Pth}_{\boldsymbol{\mathcal{B}}^{(1)}, \varphi}
\mor
[\mathrm{Pth}_{\boldsymbol{\mathcal{B}}^{(1)}}]_{\varphi}
$$
is a surjective $\Sigma^{\boldsymbol{\mathcal{A}}^{(1)}}$-homomorphism from $\mathbf{Pth}_{\boldsymbol{\mathcal{B}}^{(1)}}^{\mathbf{f}^{(1)}}$ to $[\mathbf{Pth}_{\boldsymbol{\mathcal{B}}^{(1)}}^{\mathbf{f}^{(1)}}]$.
\end{proposition}

\begin{proof}
We prove that $\mathrm{pr}_{\boldsymbol{\mathcal{B}}^{(1)}, \varphi}^{\mathrm{Ker}(\mathrm{CH}^{(1)})}$ is compatible with every operation symbol in $\Sigma^{\boldsymbol{\mathcal{A}}^{(1)}}$.

{\sffamily The mapping $\mathrm{pr}_{\boldsymbol{\mathcal{B}}^{(1)}, \varphi}^{\mathrm{Ker}(\mathrm{CH}^{(1)})}$ is a $\Sigma^{\boldsymbol{\mathcal{A}}^{(1)}}$-homomorphism.}

The fact that the mapping $\mathrm{pr}_{\boldsymbol{\mathcal{B}}^{(1)}, \varphi}^{\mathrm{Ker}(\mathrm{CH}^{(1)})}$ is a $\Sigma^{\boldsymbol{\mathcal{A}}^{(1)}}$-homomorphism follows from Remark~\ref{RDSigmaAlg}.

{\sffamily The mapping $\mathrm{pr}_{\boldsymbol{\mathcal{B}}^{(1)}, \varphi}^{\mathrm{Ker}(\mathrm{CH}^{(1)})}$ is compatible with the rewrite rules.}

Let $s$ be a sort in $S$ and $\mathfrak{P}$ a rewrite rule in $\mathcal{A}_{s}^{(1)}$. Thus,
$$
\left[
\mathfrak{p}^{\mathbf{Pth}_{\boldsymbol{\mathcal{B}}^{(1)}}^{\mathbf{f}^{(1)}}}
\right]_{\varphi(s)}
=
\left[
f_{s}^{(1)\flat}\left(
\mathfrak{p}
\right)
\right]_{\varphi(s)}
=
\mathfrak{p}^{[\mathbf{Pth}_{\boldsymbol{\mathcal{B}}^{(1)}}^{\mathbf{f}^{(1)}}]}.
$$

Hence, $\mathrm{pr}_{\boldsymbol{\mathcal{B}}^{(1)}, \varphi}^{\mathrm{Ker}(\mathrm{CH}^{(1)})}$ is compatible with the rewrite rules.

{\sffamily The mapping $\mathrm{pr}_{\boldsymbol{\mathcal{B}}^{(1)}, \varphi}^{\mathrm{Ker}(\mathrm{CH}^{(1)})}$ is compatible with the $1$-source.}

Let $s$ be a sort in $S$ and let us consider the $0$-source operation symbol $\mathrm{sc}_{s}^{0}$ in $\Sigma^{\boldsymbol{\mathcal{A}}^{(1)}}_{s,s}$. Let $\mathfrak{P}$ be a path in $\mathrm{Pth}_{\boldsymbol{\mathcal{B}}^{(1)}, \varphi(s)}$.

The following chain of equalities holds
\allowdisplaybreaks
\begin{align*}
\left[
\mathrm{sc}_{s}^{0\mathbf{Pth}_{\boldsymbol{\mathcal{B}}^{(1)}}^{\mathbf{f}^{(1)}}}\left(
\mathfrak{P}
\right)
\right]_{\varphi(s)}
&=
\left[
\mathrm{sc}_{\varphi(s)}^{0\mathbf{Pth}_{\boldsymbol{\mathcal{B}}^{(1)}}}\left(
\mathfrak{P}
\right)
\right]_{\varphi(s)}
\tag{1}
\\
&=
\mathrm{sc}_{\varphi(s)}^{0[\mathbf{Pth}_{\boldsymbol{\mathcal{B}}^{(1)}}]}\left(
\left[
\mathfrak{P}
\right]_{\varphi(s)}
\right)
\tag{2}
\\
&=
\mathrm{sc}_{s}^{0[\mathbf{Pth}_{\boldsymbol{\mathcal{B}}^{(1)}}^{\mathbf{f}^{(1)}}]}\left(
\left[
\mathfrak{P}
\right]_{\varphi(s)}
\right).
\tag{3}
\end{align*}
The first equality unravels the interpretation of the operation symbol $\mathrm{sc}_{s}^{0}$ in the partial $\Sigma^{\boldsymbol{\mathcal{A}}^{(1)}}$-algebra $\mathbf{Pth}_{\boldsymbol{\mathcal{B}}^{(1)}}^{\mathbf{f}^{(1)}}$, introduced in Proposition~\ref{PPthBCatAlg};
the second equality follows from the fact that, according to Proposition~\ref{PCHCatAlg}, $\mathrm{pr}^{\mathrm{Ker}(\mathrm{CH}^{(1)})}_{\boldsymbol{\mathcal{B}}^{(1)}}$ is a $\Lambda^{\boldsymbol{\mathcal{B}}^{(1)}}$-homomorphism from $\mathbf{Pth}_{\boldsymbol{\mathcal{B}}^{(1)}}$ to $[\mathbf{Pth}_{\boldsymbol{\mathcal{B}}^{(1)}}]$;
finally, the last equality recovers the interpretation of the operation symbol $\mathrm{sc}_{s}^{0}$ in the partial $\Sigma^{\boldsymbol{\mathcal{A}}^{(1)}}$-algebra $[\mathbf{Pth}_{\boldsymbol{\mathcal{B}}^{(1)}}^{\mathbf{f}^{(1)}}]$, introduced in Proposition~\ref{PQPthBCatAlg}.

Hence, $\mathrm{pr}_{\boldsymbol{\mathcal{B}}^{(1)}, \varphi}^{\mathrm{Ker}(\mathrm{CH}^{(1)})}$ is compatible with the $0$-source operation.

{\sffamily The mapping $\mathrm{pr}_{\boldsymbol{\mathcal{B}}^{(1)}, \varphi}^{\mathrm{Ker}(\mathrm{CH}^{(1)})}$ is compatible with the $0$-target.}

Let $s$ be a sort in $S$ and let us consider the $0$-source operation symbol $\mathrm{tg}_{s}^{0}$ in $\Sigma^{\boldsymbol{\mathcal{A}}^{(1)}}_{s,s}$. Let $\mathfrak{P}$ be a path in $\mathrm{Pth}_{\boldsymbol{\mathcal{B}}^{(1)}, \varphi(s)}$, then the following equality holds
$$
\left[
\mathrm{tg}_{s}^{0\mathbf{Pth}_{\boldsymbol{\mathcal{B}}^{(1)}}^{\mathbf{f}^{(1)}}}\left(
\mathfrak{P}
\right)
\right]_{\varphi(s)}
=
\mathrm{tg}_{s}^{0[\mathbf{Pth}_{\boldsymbol{\mathcal{B}}^{(1)}}^{\mathbf{f}^{(1)}}]}\left(
\left[
\mathfrak{P}
\right]_{\varphi(s)}
\right).
$$

The proof of this case is identical to that of the $0$-source.

Hence, $\mathrm{pr}_{\boldsymbol{\mathcal{B}}^{(1)}, \varphi}^{\mathrm{Ker}(\mathrm{CH}^{(1)})}$ is compatible with the $0$-target operation.

{\sffamily The mapping $\mathrm{pr}_{\boldsymbol{\mathcal{B}}^{(1)}, \varphi}^{\mathrm{Ker}(\mathrm{CH}^{(1)})}$ is compatible with the $0$-composition.}

Let $s$ be a sort in $S$ and let us consider the $0$-composition operation symbol $\circ_{s}^{0}$ in $\Sigma_{ss,s}^{\boldsymbol{\mathcal{A}}^{(1)}}$. Let $\mathfrak{P}$ and $\mathfrak{Q}$ be two paths in $\mathrm{Pth}_{\boldsymbol{\mathcal{B}}^{(1)}, \varphi(s)}$ such that
$$
\mathrm{sc}_{\boldsymbol{\mathcal{B}}^{(1)}, \varphi(s)}^{(0,1)}\left(\mathfrak{Q}\right)
=
\mathrm{tg}_{\boldsymbol{\mathcal{B}}^{(1)}, \varphi(s)}^{(0,1)}\left(\mathfrak{P}\right).
$$
Then the following equality holds
$$
\left[
\mathfrak{Q}
\circ_{s}^{0\mathbf{Pth}_{\boldsymbol{\mathcal{B}}^{(1)}}^{\mathbf{f}^{(1)}}}
\mathfrak{P}
\right]_{\varphi(s)}
=
\left[
\mathfrak{Q}
\right]_{\varphi(s)}
\circ_{s}^{0[\mathbf{Pth}_{\boldsymbol{\mathcal{B}}^{(1)}}^{\mathbf{f}^{(1)}}]}
\left[
\mathfrak{P}
\right]_{\varphi(s)}.
$$

The proof of this case is identical to that of the $0$-source.

Hence, $\mathrm{pr}_{\boldsymbol{\mathcal{B}}^{(1)}, \varphi}^{\mathrm{Ker}(\mathrm{CH}^{(1)})}$ is compatible with the $0$-composition operation.

This completes the proof.
\end{proof}

\section{A structure of partial $\Sigma^{\boldsymbol{\mathcal{A}}^{(1)}}$-algebra on $\mathrm{PT}_{\boldsymbol{\mathcal{B}}^{(1)}, \varphi}$}

In this section we show that, given a first-order morphism $\mathbf{f}^{(1)}$, $\mathrm{PT}_{\boldsymbol{\mathcal{B}}^{(1)}, \varphi}$ is equipped in a natural way with a structure of partial $\Sigma^{\boldsymbol{\mathcal{A}}^{(1)}}$-algebra that we denote by $\mathbf{PT}_{\boldsymbol{\mathcal{B}}^{(1)}}^{\mathbf{f}^{(1)}}$.

\begin{proposition}
\label{PPTBCatAlg}
Let $\mathbf{f}^{(1)}=(\varphi, c, (f^{(i)})_{i\in 2})$ be a first-order morphism from $\boldsymbol{\mathcal{A}}^{(1)}$ to $\boldsymbol{\mathcal{B}}^{(1)}$. Then the $S$-sorted set $\mathrm{PT}_{\boldsymbol{\mathcal{B}}^{(1)}, \varphi}$ is equipped, in a natural way, with a structure of partial $\Sigma^{\boldsymbol{\mathcal{A}}^{(1)}}$-algebra.
\end{proposition}

\begin{proof}
Let us denote by $\mathbf{PT}_{\boldsymbol{\mathcal{B}}^{(1)}}^{\mathbf{f}^{(1)}}$ the partial $\Sigma^{\boldsymbol{\mathcal{A}}^{(1)}}$-algebra defined as follows:

\textsf{(1)}
The underlying $S$-sorted set of $\mathbf{PT}_{\boldsymbol{\mathcal{B}}^{(1)}}^{\mathbf{f}^{(1)}}$ is $\mathrm{PT}_{\boldsymbol{\mathcal{B}}^{(1)}, \varphi}=(\mathrm{PT}_{\boldsymbol{\mathcal{B}}^{(1)}, \varphi(s)})_{s\in S}$.

\textsf{(2)}
For every $(\mathbf{s}, s)\in S^{\ast}\times S$ and every operation symbol $\sigma\in\Sigma_{\mathbf{s},s}$, the operation $\sigma^{\mathbf{PT}_{\boldsymbol{\mathcal{B}}^{(1)}}^{\mathbf{f}^{(1)}}}$ is given by the interpretation of $\sigma$ in the $\Sigma$-algebra $\mathbf{c}_{\mathfrak{d}}^{\ast}(\mathbf{PT}_{\boldsymbol{\mathcal{B}}^{(1)}}^{(0,1)})$. That is, its interpretation is given by the derived operation in $\mathbf{PT}^{(0,1)}_{\boldsymbol{\mathcal{B}}^{(1)}}$, the $\Lambda$-reduct of the partial $\Lambda^{\boldsymbol{\mathcal{B}}^{(1)}}$-algebra $\mathbf{PT}_{\boldsymbol{\mathcal{B}}^{(1)}}$, introduced in Proposition~\ref{PPTCatAlg}.

\textsf{(3)}
For every $s\in S$ and every $\mathfrak{p}\in\mathcal{A}_{s}^{(1)}$, the constant $\mathfrak{p}^{\mathbf{PT}_{\boldsymbol{\mathcal{B}}^{(1)}}^{\mathbf{f}^{(1)}}}$ is given by
$$
\mathfrak{p}^{\mathbf{PT}_{\boldsymbol{\mathcal{B}}^{(1)}}^{\mathbf{f}^{(1)}}}
=
\mathrm{CH}^{(1)}_{\boldsymbol{\mathcal{B}}^{(1)}, \varphi(s)}\left(
f^{(1)\flat}_{s}\left(
\mathfrak{p}^{\mathbf{Pth}_{\boldsymbol{\mathcal{A}}^{(1)}}}
\right)
\right).
$$

\textsf{(4)}
For every $s\in S$, the interpretations of the operations $\mathrm{sc}_{s}^{0}$ and $\mathrm{tg}_{s}^{0}$ are given by
\begin{align*}
\mathrm{sc}_{s}^{0\mathbf{PT}_{\boldsymbol{\mathcal{B}}^{(1)}}^{\mathbf{f}^{(1)}}}
&=
\mathrm{sc}_{\varphi(s)}^{0\mathbf{PT}_{\boldsymbol{\mathcal{B}}^{(1)}}}
&&\mbox{and}&
\mathrm{tg}_{s}^{0\mathbf{PT}_{\boldsymbol{\mathcal{B}}^{(1)}}^{\mathbf{f}^{(1)}}}
&=
\mathrm{tg}_{\varphi(s)}^{0\mathbf{PT}_{\boldsymbol{\mathcal{B}}^{(1)}}}.
\end{align*}
That is, their interpretations are given by the interpretations of $\mathrm{sc}_{\varphi(s)}^{0}$ and $\mathrm{tg}_{\varphi(s)}^{0}$ in $\mathbf{PT}_{\boldsymbol{\mathcal{B}}^{(1)}}$ that, we recall, was introduced in Proposition~\ref{PPTCatAlg}.

\textsf{(5)}
Similarly,  for every $s \in S$, the interpretations of the operation $\circ_{s}^{0}$ is given by
$$
\circ_{s}^{0\mathbf{PT}_{\boldsymbol{\mathcal{B}}^{(1)}}^{\mathbf{f}^{(1)}}}
=
\circ_{\varphi(s)}^{0\mathbf{PT}_{\boldsymbol{\mathcal{B}}^{(1)}}}.
$$
That is, its interpretation is given by the interpretation of $\circ_{\varphi(s)}^{0}$ in $\mathbf{PT}_{\boldsymbol{\mathcal{B}}^{(1)}}$ that, we recall, was introduced in Proposition~\ref{PPTCatAlg}.
\end{proof}

Now, given a first-order morphism $\mathbf{f}^{(1)}$, $[\mathrm{PT}_{\boldsymbol{\mathcal{B}}^{(1)}}]_{\varphi}$ is equipped in a natural way with a structure of partial $\Sigma^{\boldsymbol{\mathcal{A}}^{(1)}}$-algebra that we denote by $[\mathbf{PT}_{\boldsymbol{\mathcal{B}}^{(1)}}^{\mathbf{f}^{(1)}}]$. Thus, we show that $\mathrm{pr}^{\Theta^{[1]}}_{\boldsymbol{\mathcal{B}}^{(1)}, \varphi}$ is a surjective $\Sigma^{\boldsymbol{\mathcal{A}}^{(1)}}$-homomorphism from $\mathbf{PT}_{\boldsymbol{\mathcal{B}}^{(1)}}^{\mathbf{f}^{(1)}}$ to $[\mathbf{PT}_{\boldsymbol{\mathcal{B}}^{(1)}}^{\mathbf{f}^{(1)}}]$.

\begin{proposition}
\label{PQPTBCatAlg}
Let $\mathbf{f}^{(1)}=(\varphi, c, (f^{(i)})_{i\in 2})$ be a first-order morphism from $\boldsymbol{\mathcal{A}}^{(1)}$ to $\boldsymbol{\mathcal{B}}^{(1)}$. Then the $S$-sorted set $[\mathrm{PT}_{\boldsymbol{\mathcal{B}}^{(1)}}]_{\varphi}$ is equipped, in a natural way, with a structure of partial $\Sigma^{\boldsymbol{\mathcal{A}}^{(1)}}$-algebra.
\end{proposition}

\begin{proof}
Let us denote by $[\mathbf{PT}_{\boldsymbol{\mathcal{B}}^{(1)}}^{\mathbf{f}^{(1)}}]$ the $\Sigma^{\boldsymbol{\mathcal{A}}^{(1)}}$-algebra defined as follows:

\textsf{(1)}
The underlying $S$-sorted set of $[\mathbf{PT}_{\boldsymbol{\mathcal{B}}^{(1)}}^{\mathbf{f}^{(1)}}]$ is $[\mathrm{PT}_{\boldsymbol{\mathcal{B}}^{(1)}}]_{\varphi}=([\mathrm{PT}_{\boldsymbol{\mathcal{B}}^{(1)}}]_{\varphi(s)})_{s\in S}$.

\textsf{(2)}
For every $(\mathbf{s}, s)\in S^{\ast}\times S$ and every operation symbol $\sigma\in\Sigma_{\mathbf{s},s}$, the operation $\sigma^{\mathbf{PT}_{\boldsymbol{\mathcal{B}}^{(1)}}^{\mathbf{f}^{(1)}}}$ is given by the interpretation of $\sigma$ in the $\Sigma$-algebra $\mathbf{c}_{\mathfrak{d}}^{\ast}([\mathbf{PT}_{\boldsymbol{\mathcal{B}}^{(1)}}^{(0,1)}])$. That is, its interpretation is given by the derived operation in $[\mathbf{PT}^{(0,1)}_{\boldsymbol{\mathcal{B}}^{(1)}}]$, the $\Lambda$-reduct of the partial $\Lambda^{\boldsymbol{\mathcal{B}}^{(1)}}$-algebra $[\mathbf{PT}_{\boldsymbol{\mathcal{B}}^{(1)}}]$, introduced in Proposition~\ref{PPTQCatAlg}.

\textsf{(3)}
For every $s\in S$ and every $\mathfrak{p}\in\mathcal{A}_{s}^{(1)}$, the constant $\mathfrak{p}^{\mathbf{PT}_{\boldsymbol{\mathcal{B}}^{(1)}}^{\mathbf{f}^{(1)}}}$ is given by
$$
\mathfrak{p}^{[\mathbf{PT}_{\boldsymbol{\mathcal{B}}^{(1)}}^{\mathbf{f}^{(1)}}]}
=
\left[
\mathrm{CH}^{(1)}_{\boldsymbol{\mathcal{B}}^{(1)}, \varphi(s)}\left(
f^{(1)\flat}_{s}\left(
\mathfrak{p}^{\mathbf{Pth}_{\boldsymbol{\mathcal{A}}^{(1)}}}
\right)
\right)
\right]_{\varphi(s)}.
$$

\textsf{(4)}
For every $s\in S$, the interpretations of the operations $\mathrm{sc}_{s}^{0}$ and $\mathrm{tg}_{s}^{0}$ are given by
\begin{align*}
\mathrm{sc}_{s}^{0[\mathbf{PT}_{\boldsymbol{\mathcal{B}}^{(1)}}^{\mathbf{f}^{(1)}}]}
&=
\mathrm{sc}_{\varphi(s)}^{0[\mathbf{PT}_{\boldsymbol{\mathcal{B}}^{(1)}}]}
&&\mbox{and}&
\mathrm{tg}_{s}^{0[\mathbf{PT}_{\boldsymbol{\mathcal{B}}^{(1)}}^{\mathbf{f}^{(1)}}]}
&=
\mathrm{tg}_{\varphi(s)}^{0[\mathbf{PT}_{\boldsymbol{\mathcal{B}}^{(1)}}]}.
\end{align*}
That is, their interpretations are given by the interpretations of $\mathrm{sc}_{\varphi(s)}^{0}$ and $\mathrm{tg}_{\varphi(s)}^{0}$ in $[\mathbf{PT}_{\boldsymbol{\mathcal{B}}^{(1)}}]$ that, we recall, was introduced in Proposition~\ref{PPTQCatAlg}.

\textsf{(5)}
Similarly,  for every $s \in S$, the interpretations of the operation $\circ_{s}^{0}$ is given by
$$
\circ_{s}^{0[\mathbf{PT}_{\boldsymbol{\mathcal{B}}^{(1)}}^{\mathbf{f}^{(1)}}]}
=
\circ_{\varphi(s)}^{0[\mathbf{PT}_{\boldsymbol{\mathcal{B}}^{(1)}}]}.
$$
That is, its interpretation is given by the interpretation of $\circ_{\varphi(s)}^{0}$ in $[\mathbf{PT}_{\boldsymbol{\mathcal{B}}^{(1)}}]$ that, we recall, was introduced in Proposition~\ref{PPTQCatAlg}.
\end{proof}

\begin{proposition}
\label{PPrPTBCatHom}
Let $\mathbf{f}^{(1)}=(\varphi, c, (f^{(i)})_{i\in 2})$ be a first-order morphism from $\boldsymbol{\mathcal{A}}^{(1)}$ to $\boldsymbol{\mathcal{B}}^{(1)}$. Then the mapping 
$$
\mathrm{pr}_{\boldsymbol{\mathcal{B}}^{(1)}, \varphi}^{\Theta^{[1]}}
\colon 
\mathrm{PT}_{\boldsymbol{\mathcal{B}}^{(1)}, \varphi}
\mor
[\mathrm{PT}_{\boldsymbol{\mathcal{B}}^{(1)}}]_{\varphi}
$$
is a surjective $\Sigma^{\boldsymbol{\mathcal{A}}^{(1)}}$-homomorphism from $\mathbf{PT}_{\boldsymbol{\mathcal{B}}^{(1)}}^{\mathbf{f}^{(1)}}$ to $[\mathbf{PT}_{\boldsymbol{\mathcal{B}}^{(1)}}^{\mathbf{f}^{(1)}}]$.
\end{proposition}

\begin{proof}
We prove that $\mathrm{pr}_{\boldsymbol{\mathcal{B}}^{(1)}, \varphi}^{\Theta^{[1]}}$ is compatible with every operation symbol in $\Sigma^{\boldsymbol{\mathcal{A}}^{(1)}}$.

{\sffamily The mapping $\mathrm{pr}_{\boldsymbol{\mathcal{B}}^{(1)}, \varphi}^{\Theta^{[1]}}$ is a $\Sigma$-homomorphism.}

Note that $\mathrm{pr}_{\boldsymbol{\mathcal{B}}^{(1)}, \varphi}^{\Theta^{[1]}} = \mathbf{c}_{\mathfrak{d}}^{\ast} (\mathrm{pr}_{\boldsymbol{\mathcal{B}}^{(1)}}^{\Theta^{[1]}})$. By Proposition~\ref{CPTQPr}, the mapping $\mathrm{pr}_{\boldsymbol{\mathcal{B}}^{(1)}}^{\Theta^{[1]}}$ is a $\Lambda^{\boldsymbol{\mathcal{B}}^{(1)}}$-homomorphism, thus in particular a $\Lambda$-homomorphism. Therefore, it follows from Proposition~\ref{PFunSig} that the mapping $\mathrm{pr}_{\boldsymbol{\mathcal{B}}^{(1)}, \varphi}^{\Theta^{[1]}}$ is a $\Sigma$-homomorphism.


{\sffamily The mapping $\mathrm{pr}_{\boldsymbol{\mathcal{B}}^{(1)}, \varphi}^{\Theta^{[1]}}$ is compatible with the rewrite rules.}

Let $s$ be a sort in $S$ and $\mathfrak{p}$ a rewrite rule in $\mathcal{A}_{s}^{(1)}$. Thus,
$$
\left[
\mathfrak{p}^{\mathbf{PT}_{\boldsymbol{\mathcal{B}}^{(1)}}^{\mathbf{f}^{(1)}}}
\right]_{\varphi(s)}
=
\left[
\mathrm{CH}^{(1)}_{\boldsymbol{\mathcal{B}}^{(1)}, \varphi(s)}\left(
f^{(1)\flat}_{s}\left(
\mathfrak{p}^{\mathbf{Pth}_{\boldsymbol{\mathcal{A}}^{(1)}}}
\right)
\right)
\right]_{\varphi(s)}
=
\mathfrak{p}^{[\mathbf{PT}_{\boldsymbol{\mathcal{B}}^{(1)}}^{\mathbf{f}^{(1)}}]}.
$$

Hence, $\mathrm{pr}_{\boldsymbol{\mathcal{B}}^{(1)}, \varphi}^{\Theta^{[1]}}$ is compatible with the rewrite rules.

{\sffamily The mapping $\mathrm{pr}_{\boldsymbol{\mathcal{B}}^{(1)}, \varphi}^{\Theta^{[1]}}$ is compatible with the $0$-source.}

Let $s$ be a sort in $S$ and let us consider the $0$-source operation symbol $\mathrm{sc}_{s}^{0}$ in $\Sigma^{\boldsymbol{\mathcal{A}}^{(1)}}_{s,s}$. Let $P$ be a path term in $\mathrm{PT}_{\boldsymbol{\mathcal{B}}^{(1)}, \varphi(s)}$.

The following chain of equalities holds
\allowdisplaybreaks
\begin{align*}
\left[
\mathrm{sc}_{s}^{0\mathbf{PT}_{\boldsymbol{\mathcal{B}}^{(1)}}^{\mathbf{f}^{(1)}}}\left(
P
\right)
\right]_{\varphi(s)}
&=
\left[
\mathrm{sc}_{\varphi(s)}^{0\mathbf{PT}_{\boldsymbol{\mathcal{B}}^{(1)}}}\left(
P
\right)
\right]_{\varphi(s)}
\tag{1}
\\
&=
\mathrm{sc}_{\varphi(s)}^{0[\mathbf{Pth}_{\boldsymbol{\mathcal{B}}^{(1)}}]}\left(
\left[
P
\right]_{\varphi(s)}
\right)
\tag{2}
\\
&=
\mathrm{sc}_{s}^{0[\mathbf{Pth}_{\boldsymbol{\mathcal{B}}^{(1)}}^{\mathbf{f}^{(1)}}]}\left(
\left[
P
\right]_{\varphi(s)}
\right).
\tag{3}
\end{align*}
The first equality unravels the interpretation of the operation symbol $\mathrm{sc}_{s}^{0}$ in the partial $\Sigma^{\boldsymbol{\mathcal{A}}^{(1)}}$-algebra $\mathbf{PT}_{\boldsymbol{\mathcal{B}}^{(1)}}^{\mathbf{f}^{(1)}}$, introduced in Proposition~\ref{PPTBCatAlg};
the second equality follows from the fact that, according to Proposition~\ref{CPTQPr}, $\mathrm{pr}^{\Theta^{[1]}}_{\boldsymbol{\mathcal{B}}^{(1)}}$ is a $\Lambda^{\boldsymbol{\mathcal{B}}^{(1)}}$-homomorphism from $\mathbf{PT}_{\boldsymbol{\mathcal{B}}^{(1)}}$ to $[\mathbf{PT}_{\boldsymbol{\mathcal{B}}^{(1)}}]$;
finally, the last equality recovers the interpretation of the operation symbol $\mathrm{sc}_{s}^{0}$ in the partial $\Sigma^{\boldsymbol{\mathcal{A}}^{(1)}}$-algebra $[\mathbf{PT}_{\boldsymbol{\mathcal{B}}^{(1)}}^{\mathbf{f}^{(1)}}]$, introduced in Proposition~\ref{PQPTBCatAlg}.

Hence, $\mathrm{pr}_{\boldsymbol{\mathcal{B}}^{(1)}, \varphi}^{\Theta^{[1]}}$ is compatible with the $0$-source operation.

{\sffamily The mapping $\mathrm{pr}_{\boldsymbol{\mathcal{B}}^{(1)}, \varphi}^{\Theta^{[1]}}$ is compatible with the $0$-target.}

Let $s$ be a sort in $S$ and let us consider the $0$-source operation symbol $\mathrm{tg}_{s}^{0}$ in $\Sigma^{\boldsymbol{\mathcal{A}}^{(1)}}_{s,s}$. Let $P$ be a path term in $\mathrm{PT}_{\boldsymbol{\mathcal{B}}^{(1)}, \varphi(s)}$, then the following equality holds
$$
\left[
\mathrm{tg}_{s}^{0\mathbf{PT}_{\boldsymbol{\mathcal{B}}^{(1)}}^{\mathbf{f}^{(1)}}}\left(
P
\right)
\right]_{\varphi(s)}
=
\mathrm{tg}_{s}^{0[\mathbf{PT}_{\boldsymbol{\mathcal{B}}^{(1)}}^{\mathbf{f}^{(1)}}]}\left(
\left[
P
\right]_{\varphi(s)}
\right).
$$

The proof of this case is identical to that of the $0$-source.

Hence, $\mathrm{pr}_{\boldsymbol{\mathcal{B}}^{(1)}, \varphi}^{\Theta^{[1]}}$ is compatible with the $0$-target operation.

{\sffamily The mapping $\mathrm{pr}_{\boldsymbol{\mathcal{B}}^{(1)}, \varphi}^{\Theta^{[1]}}$ is compatible with the $0$-composition.}

Let $s$ be a sort in $S$ and let us consider the $0$-composition operation symbol $\circ_{s}^{0}$ in $\Sigma_{ss,s}^{\boldsymbol{\mathcal{A}}^{(1)}}$. Let $P$ and $Q$ be two path terms in $\mathrm{PT}_{\boldsymbol{\mathcal{B}}^{(1)}, \varphi(s)}$ such that
$$
\mathrm{sc}_{\boldsymbol{\mathcal{B}}^{(1)}, \varphi(s)}^{(0,1)}\left(
\mathrm{ip}^{(1,Y)@}_{\boldsymbol{\mathcal{B}}^{(1)}, \varphi(s)}\left(
Q
\right)
\right)
=
\mathrm{tg}_{\boldsymbol{\mathcal{B}}^{(1)}, \varphi(s)}^{(0,1)}\left(
\mathrm{ip}^{(1,Y)@}_{\boldsymbol{\mathcal{B}}^{(1)}, \varphi(s)}\left(
P
\right)
\right).
$$
Then the following equality holds
$$
\left[
Q
\circ_{s}^{0\mathbf{Pth}_{\boldsymbol{\mathcal{B}}^{(1)}}^{\mathbf{f}^{(1)}}}
P
\right]_{\varphi(s)}
=
\left[
Q
\right]_{\varphi(s)}
\circ_{s}^{0[\mathbf{Pth}_{\boldsymbol{\mathcal{B}}^{(1)}}^{\mathbf{f}^{(1)}}]}
\left[
P
\right]_{\varphi(s)}.
$$

The proof of this case is identical to that of the $0$-source.

Hence, $\mathrm{pr}_{\boldsymbol{\mathcal{B}}^{(1)}, \varphi}^{\Theta^{[1]}}$ is compatible with the $0$-composition operation.

This completes the proof.
\end{proof}

Finally, we show that the partial $\Sigma^{\boldsymbol{\mathcal{A}}^{(1)}}$-algebras $[\mathbf{Pth}_{\boldsymbol{\mathcal{B}}^{(1)}}^{\mathbf{f}^{(1)}}]$ and $[\mathbf{PT}_{\boldsymbol{\mathcal{B}}^{(1)}}^{\mathbf{f}^{(1)}}]$ are isomorphic. Indeed the mappings $\mathrm{ip}^{([1], Y)@}_{\boldsymbol{\mathcal{B}}^{(1)}, \varphi}$ and $\mathrm{CH}^{[ 1 ]}_{\boldsymbol{\mathcal{B}}^{(1)}, \varphi}$ are a pair of inverse $\Sigma^{\boldsymbol{\mathcal{A}}^{(1)}}$-isomorphisms.

\begin{proposition}
\label{PQCHBCatHom}
Let $\mathbf{f}^{(1)}=(\varphi, c, (f^{(i)})_{i\in 2})$ be a first-order morphism from $\boldsymbol{\mathcal{A}}^{(1)}$ to $\boldsymbol{\mathcal{B}}^{(1)}$. Then the mapping 
$$
\mathrm{CH}^{[1]}_{\boldsymbol{\mathcal{B}}^{(1)}, \varphi}
\colon 
[\mathrm{Pth}_{\boldsymbol{\mathcal{B}}^{(1)}}]_{\varphi}
\mor
[\mathrm{PT}_{\boldsymbol{\mathcal{B}}^{(1)}}]_{\varphi}
$$
is a $\Sigma^{\boldsymbol{\mathcal{A}}^{(1)}}$-homomorphism from $[\mathbf{Pth}_{\boldsymbol{\mathcal{B}}^{(1)}}^{\mathbf{f}^{(1)}}]$ to $[\mathbf{PT}_{\boldsymbol{\mathcal{B}}^{(1)}}^{\mathbf{f}^{(1)}}]$.
\end{proposition}

\begin{proof}
We prove that $\mathrm{CH}^{[1]}_{\boldsymbol{\mathcal{B}}^{(1)}, \varphi}$ is compatible with every operation symbol in $\Sigma^{\boldsymbol{\mathcal{A}}^{(1)}}$.

{\sffamily The mapping $\mathrm{CH}^{[1]}_{\boldsymbol{\mathcal{B}}^{(1)}, \varphi}$ is a $\Sigma$-homomorphism.}

Note that $\mathrm{CH}^{[1]}_{\boldsymbol{\mathcal{B}}^{(1)}, \varphi} = \mathbf{c}_{\mathfrak{d}}^{\ast} (\mathrm{CH}^{[1]}_{\boldsymbol{\mathcal{B}}^{(1)}})$. By Proposition~\ref{CIso}, the mapping $\mathrm{CH}^{[1]}_{\boldsymbol{\mathcal{B}}^{(1)}}$ is a $\Lambda^{\boldsymbol{\mathcal{B}}^{(1)}}$-homomorphism, thus in particular a $\Lambda$-homomorphism. Therefore, it follows from Proposition~\ref{PFunSig} that the mapping $\mathrm{CH}^{[1]}_{\boldsymbol{\mathcal{B}}^{(1)}, \varphi}$ is a $\Sigma$-homomorphism.

{\sffamily The mapping $\mathrm{CH}^{[1]}_{\boldsymbol{\mathcal{B}}^{(1)}, \varphi}$ is compatible with the rewrite rules.}

Let $s$ be a sort in $S$ and $\mathfrak{p}$ a rewrite rule in $\mathcal{A}_{s}^{(1)}$. Thus,
\begin{align*}
\mathrm{CH}^{[1]}_{\boldsymbol{\mathcal{B}}^{(1)}, \varphi(s)}\left(
\mathfrak{p}^{[\mathbf{Pth}_{\boldsymbol{\mathcal{B}}^{(1)}}^{\mathbf{f}^{(1)}}]}
\right)
&=
\mathrm{CH}^{[1]}_{\boldsymbol{\mathcal{B}}^{(1)}, \varphi(s)}\left(
\left[
f^{(1)\flat}\left(
\mathrm{ech}^{(1,\mathcal{A}^{(1)})}_{\boldsymbol{\mathcal{A}}^{(1)}, s} \left(
\mathfrak{p}
\right)
\right)
\right]_{\varphi(s)}
\right)
\tag{1}
\\
&=
\left[
\mathrm{CH}^{(1)}_{\boldsymbol{\mathcal{B}}^{(1)}, \varphi(s)}\left(
f^{(1)\flat}\left(
\mathrm{ech}^{(1,\mathcal{A}^{(1)})}_{\boldsymbol{\mathcal{A}}^{(1)}, s} \left(
\mathfrak{p}
\right)
\right)
\right)
\right]_{\varphi(s)}
\tag{2}
\\
&=
\mathfrak{p}^{[\mathbf{PT}_{\boldsymbol{\mathcal{B}}^{(1)}}^{\mathbf{f}^{(1)}}]}
\tag{3}
\end{align*}
The first equality unravels the definition of the constant $\mathfrak{p}$ in the partial $\Sigma^{\boldsymbol{\mathcal{A}}^{(1)}}$-algebra $[\mathbf{Pth}_{\boldsymbol{\mathcal{B}}^{(1)}}^{\mathbf{f}^{(1)}}]$, introduced in Proposition~\ref{PQPthBCatAlg};
the second equality unravels the definition of the mapping $\mathrm{CH}^{[1]}_{\boldsymbol{\mathcal{B}}^{(1)}}$, introduced in Definition~\ref{DPTQCH};
finally, the last equality recovers the definition of the constant $\mathfrak{p}$ in the partial $\Sigma^{\boldsymbol{\mathcal{A}}^{(1)}}$-algebra $[\mathbf{PT}_{\boldsymbol{\mathcal{B}}^{(1)}}^{\mathbf{f}^{(1)}}]$, introduced in Proposition~\ref{PQPTBCatAlg};

Hence, $\mathrm{CH}^{[1]}_{\boldsymbol{\mathcal{B}}^{(1)}, \varphi}$ is compatible with the rewrite rules.

{\sffamily The mapping $\mathrm{CH}^{[1]}_{\boldsymbol{\mathcal{B}}^{(1)}, \varphi}$ is compatible with the $0$-source.}

Let $s$ be a sort in $S$ and let us consider the $0$-source operation symbol $\mathrm{sc}_{s}^{0}$ in $\Sigma^{\boldsymbol{\mathcal{A}}^{(1)}}_{s,s}$. Let $\mathfrak{P}$ be a path in $\mathrm{Pth}_{\boldsymbol{\mathcal{B}}^{(1)}, \varphi(s)}$.

The following chain of equalities holds
\allowdisplaybreaks
\begin{align*}
\mathrm{CH}^{[1]}_{\boldsymbol{\mathcal{B}}^{(1)}, \varphi(s)} \left(
\mathrm{sc}_{s}^{0[\mathbf{Pth}_{\boldsymbol{\mathcal{B}}^{(1)}}^{\mathbf{f}^{(1)}}]}\left(
\left[
\mathfrak{P}
\right]_{\varphi(s)}
\right)
\right)
&=
\mathrm{CH}^{[1]}_{\boldsymbol{\mathcal{B}}^{(1)}, \varphi(s)} \left(
\mathrm{sc}_{\varphi(s)}^{0[\mathbf{Pth}_{\boldsymbol{\mathcal{B}}^{(1)}}]}\left(
\left[
\mathfrak{P}
\right]_{\varphi(s)}
\right)
\right)
\tag{1}
\\
&=
\mathrm{sc}_{\varphi(s)}^{0[\mathbf{PT}_{\boldsymbol{\mathcal{B}}^{(1)}}]}\left(
\mathrm{CH}^{[1]}_{\boldsymbol{\mathcal{B}}^{(1)}, \varphi(s)} \left(
\left[
\mathfrak{P}
\right]_{\varphi(s)}
\right)
\right)
\tag{2}
\\
&=
\mathrm{sc}_{s}^{0[\mathbf{PT}_{\boldsymbol{\mathcal{B}}^{(1)}}^{\mathbf{f}^{(1)}}]}\left(
\mathrm{CH}^{[1]}_{\boldsymbol{\mathcal{B}}^{(1)}, \varphi(s)} \left(
\left[
\mathfrak{P}
\right]_{\varphi(s)}
\right)
\right).
\tag{3}
\end{align*}

The first equality unravels the interpretation of the operation symbol $\mathrm{sc}_{s}^{0}$ in the partial $\Sigma^{\boldsymbol{\mathcal{A}}^{(1)}}$-algebra $[\mathbf{Pth}_{\boldsymbol{\mathcal{B}}^{(1)}}^{\mathbf{f}^{(1)}}]$, introduced in Proposition~\ref{PPthBCatAlg};
the second equality follows from the fact that, according to Claim~\ref{CIso}, $\mathrm{CH}^{[1]}_{\boldsymbol{\mathcal{B}}^{(1)}, \varphi}$ is a $\Lambda^{\boldsymbol{\mathcal{B}}^{(1)}}$-homomorphism from $[\mathbf{Pth}_{\boldsymbol{\mathcal{B}}^{(1)}}]$ to $[\mathbf{PT}_{\boldsymbol{\mathcal{B}}^{(1)}}]$;
finally, the last equality recovers the interpretation of the operation symbol $\mathrm{sc}_{s}^{0}$ in the partial $\Sigma^{\boldsymbol{\mathcal{A}}^{(1)}}$-algebra $[\mathbf{PT}_{\boldsymbol{\mathcal{B}}^{(1)}}^{\mathbf{f}^{(1)}}]$, introduced in Proposition~\ref{PQPTBCatAlg}.

Hence, $\mathrm{CH}^{[1]}_{\boldsymbol{\mathcal{B}}^{(1)}, \varphi}$ is compatible with the $0$-source operation.

{\sffamily The mapping $\mathrm{CH}^{[1]}_{\boldsymbol{\mathcal{B}}^{(1)}, \varphi}$ is compatible with the $0$-target.}

Let $s$ be a sort in $S$ and let us consider the $0$-target operation symbol $\mathrm{tg}_{s}^{0}$ in $\Sigma^{\boldsymbol{\mathcal{A}}^{(1)}}_{s,s}$. Let $\mathfrak{P}$ be a path in $\mathrm{Pth}_{\boldsymbol{\mathcal{B}}^{(1)}, \varphi(s)}$, then the following equality holds
\begin{multline*}
\mathrm{CH}^{[1]}_{\boldsymbol{\mathcal{B}}^{(1)}, \varphi(s)}\left(
\mathrm{tg}_{s}^{0[\mathbf{Pth}_{\boldsymbol{\mathcal{B}}^{(1)}}^{\mathbf{f}^{(1)}}]}\left(
\left[
\mathfrak{P}
\right]_{\varphi(s)}
\right)
\right)
=
\mathrm{tg}_{s}^{0[\mathbf{PT}_{\boldsymbol{\mathcal{B}}^{(1)}}^{\mathbf{f}^{(1)}}]}\left(
\mathrm{CH}^{[1]}_{\boldsymbol{\mathcal{B}}^{(1)}, \varphi(s)}\left(
\left[
\mathfrak{P}
\right]_{\varphi(s)}
\right)
\right).
\end{multline*}

The proof of this case is identical to that of the $0$-source.

Hence, $\mathrm{CH}^{[1]}_{\boldsymbol{\mathcal{B}}^{(1)}, \varphi}$ is compatible with the $0$-target operation.

{\sffamily The mapping $\mathrm{CH}^{[1]}_{\boldsymbol{\mathcal{B}}^{(1)}, \varphi}$ is compatible with the $0$-composition.}

Let $s$ be a sort in $S$ and let us consider the $0$-composition operation symbol $\circ_{s}^{0}$ in $\Sigma_{ss,s}^{\boldsymbol{\mathcal{A}}^{(1)}}$. Let $\mathfrak{P}$ and $\mathfrak{Q}$ be two paths in $\mathrm{Pth}_{\boldsymbol{\mathcal{B}}^{(1)}, \varphi(s)}$ such that
$$
\mathrm{sc}_{\boldsymbol{\mathcal{B}}^{(1)}, \varphi(s)}^{(0,2)}\left(\mathfrak{Q}\right)
=
\mathrm{tg}_{\boldsymbol{\mathcal{B}}^{(1)}, \varphi(s)}^{(0,2)}\left(\mathfrak{P}\right).
$$

Then the following equality holds
\begin{multline*}
\mathrm{CH}^{[1]}_{\boldsymbol{\mathcal{B}}^{(1)}, \varphi(s)}\left(
\left[
\mathfrak{Q}
\right]_{\varphi(s)}
\circ_{s}^{0[\mathbf{Pth}_{\boldsymbol{\mathcal{B}}^{(1)}}^{\mathbf{f}^{(1)}}]}
\left[
\mathfrak{P}
\right]_{\varphi(s)}
\right)
\\
=
\mathrm{CH}^{[1]}_{\boldsymbol{\mathcal{B}}^{(1)}, \varphi(s)}\left(
\left[
\mathfrak{Q}
\right]_{\varphi(s)}
\right)
\circ_{s}^{0[\mathbf{PT}_{\boldsymbol{\mathcal{B}}^{(1)}}^{\mathbf{f}^{(1)}}]}
\mathrm{CH}^{[1]}_{\boldsymbol{\mathcal{B}}^{(1)}, \varphi(s)}\left(
\left[
\mathfrak{P}
\right]_{\varphi(s)}
\right).
\end{multline*}

The proof of this case is identical to that of the $0$-source.

Hence, $\mathrm{CH}^{[1]}_{\boldsymbol{\mathcal{B}}^{(1)}, \varphi}$ is compatible with the $0$-composition operation.

This completes the proof.
\end{proof}

\begin{proposition}
\label{PQIpfcBCatHom}
Let $\mathbf{f}^{(1)}=(\varphi, c, (f^{(i)})_{i\in 2})$ be a first-order morphism from $\boldsymbol{\mathcal{A}}^{(1)}$ to $\boldsymbol{\mathcal{B}}^{(1)}$. Then the mapping 
$$
\mathrm{ip}^{([1], Y)@}_{\boldsymbol{\mathcal{B}}^{(1)}, \varphi}
\colon 
[\mathrm{PT}_{\boldsymbol{\mathcal{B}}^{(1)}}]_{\varphi}
\mor
[\mathrm{Pth}_{\boldsymbol{\mathcal{B}}^{(1)}}]_{\varphi}
$$
is a $\Sigma^{\boldsymbol{\mathcal{A}}^{(1)}}$-homomorphism from $[\mathbf{PT}_{\boldsymbol{\mathcal{B}}^{(1)}}^{\mathbf{f}^{(1)}}]$ to $[\mathbf{Pth}_{\boldsymbol{\mathcal{B}}^{(1)}}^{\mathbf{f}^{(1)}}]$.
\end{proposition}

\begin{proof}
We prove that $\mathrm{ip}^{([1], Y)@}_{\boldsymbol{\mathcal{B}}^{(1)}, \varphi}$ is compatible with every operation symbol in $\Sigma^{\boldsymbol{\mathcal{A}}^{(1)}}$.

{\sffamily The mapping $\mathrm{ip}^{([1], Y)@}_{\boldsymbol{\mathcal{B}}^{(1)}, \varphi}$ is a $\Sigma$-homomorphism.}

Note that $\mathrm{ip}^{([1], Y)@}_{\boldsymbol{\mathcal{B}}^{(1)}, \varphi} = \mathbf{c}_{\mathfrak{d}}^{\ast} (\mathrm{ip}^{([1], Y)@}_{\boldsymbol{\mathcal{B}}^{(1)}})$. By Proposition~\ref{CIsoIpfc}, the mapping $\mathrm{ip}^{([1], Y)@}_{\boldsymbol{\mathcal{B}}^{(1)}}$ is a $\Lambda^{\boldsymbol{\mathcal{B}}^{(1)}}$-homomorphism, thus in particular a $\Lambda$-homomorphism. Therefore, it follows from Proposition~\ref{PFunSig} that the mapping $\mathrm{ip}^{([1], Y)@}_{\boldsymbol{\mathcal{B}}^{(1)}, \varphi}$ is a $\Sigma$-homomorphism.

{\sffamily The mapping $\mathrm{ip}^{([1], X)@}_{\boldsymbol{\mathcal{B}}^{(1)}, \varphi}$ is compatible with the rewrite rules.}

Let $s$ be a sort in $S$ and $\mathfrak{p}$ a rewrite rule in $\mathcal{A}_{s}^{(1)}$. Thus,
\begin{flushleft}
$
\mathrm{ip}^{([1], Y)@}_{\boldsymbol{\mathcal{B}}^{(1)}, \varphi(s)}\left(
\mathfrak{p}^{[\mathbf{PT}_{\boldsymbol{\mathcal{B}}^{(1)}}^{\mathbf{f}^{(1)}}]}
\right)
$
\allowdisplaybreaks
\begin{align*}
&=
\mathrm{ip}^{([1], Y)@}_{\boldsymbol{\mathcal{B}}^{(1)}, \varphi(s)}\left(
\left[
\mathrm{CH}^{(1)}_{\boldsymbol{\mathcal{B}}^{(1)}, \varphi(s)}\left(
f^{(1)\flat}_{s}\left(
\mathrm{ech}^{(1,\mathcal{A}^{(1)})}_{\boldsymbol{\mathcal{A}}^{(1)}, s}\left(
\mathfrak{p}
\right)
\right)
\right)
\right]_{\varphi(s)}
\right)
\tag{1}
\\
&=
\left[
\mathrm{ip}^{(1, Y)@}_{\boldsymbol{\mathcal{B}}^{(1)}, \varphi(s)}\left(
\mathrm{CH}^{(1)}_{\boldsymbol{\mathcal{B}}^{(1)}, \varphi(s)}\left(
f^{(1)\flat}_{s}\left(
\mathrm{ech}^{(1,\mathcal{A}^{(1)})}_{\boldsymbol{\mathcal{A}}^{(1)}, s}\left(
\mathfrak{p}
\right)
\right)
\right)
\right)
\right]_{\varphi(s)}
\tag{2}
\\
&=
\left[
f^{(1)\flat}_{s}\left(
\mathrm{ech}^{(1,\mathcal{A}^{(1)})}_{\boldsymbol{\mathcal{A}}^{(1)}, s}\left(
\mathfrak{p}
\right)
\right)
\right]_{\varphi(s)}
\tag{3}
\\
&=
\mathfrak{p}^{[\mathbf{Pth}_{\boldsymbol{\mathcal{B}}^{(1)}}^{\mathbf{f}^{(1)}}]}
\tag{4}
\end{align*}
\end{flushleft}

The first equality unravels the definition of the constant $\mathfrak{p}$ in the partial $\Sigma^{\boldsymbol{\mathcal{A}}^{(1)}}$-algebra $[\mathbf{PT}_{\boldsymbol{\mathcal{B}}^{(1)}}^{\mathbf{f}^{(1)}}]$, introduced in Proposition~\ref{PQPTBCatAlg};
the second equality unravels the definition of the mapping $\mathrm{ip}^{([1], Y)@}_{\boldsymbol{\mathcal{B}}^{(1)}}$, introduced in Definition~\ref{DPTQIp};
the third equality follows from Proposition~\ref{PIpCH};
finally, the last equality recovers the definition of the constant $\mathfrak{p}$ in the partial $\Sigma^{\boldsymbol{\mathcal{A}}^{(1)}}$-algebra $[\mathbf{Pth}_{\boldsymbol{\mathcal{B}}^{(1)}}^{\mathbf{f}^{(1)}}]$, introduced in Proposition~\ref{PQPthBCatAlg}.

Hence, $\mathrm{ip}^{([1], Y)@}_{\boldsymbol{\mathcal{B}}^{(1)}, \varphi}$ is compatible with the rewrite rules.

{\sffamily The mapping $\mathrm{ip}^{([1], Y)@}_{\boldsymbol{\mathcal{B}}^{(1)}, \varphi}$ is compatible with the $0$-source.}

Let $s$ be a sort in $S$ and let us consider the $0$-source operation symbol $\mathrm{sc}_{s}^{0}$ in $\Sigma^{\boldsymbol{\mathcal{A}}^{(1)}}_{s,s}$. Let $P$ be a path term in $\mathrm{PT}_{\boldsymbol{\mathcal{B}}^{(1)}, \varphi(s)}$.

\allowdisplaybreaks
\begin{align*}
\mathrm{ip}^{([1], Y)@}_{\boldsymbol{\mathcal{B}}^{(1)}, \varphi(s)} \left(
\mathrm{sc}_{s}^{0[\mathbf{PT}_{\boldsymbol{\mathcal{B}}^{(1)}}^{\mathbf{f}^{(1)}}]}\left(
\left[
P
\right]_{\varphi(s)}
\right)
\right)
&=
\mathrm{ip}^{([1], Y)@}_{\boldsymbol{\mathcal{B}}^{(1)}, \varphi(s)} \left(
\mathrm{sc}_{\varphi(s)}^{0[\mathbf{PT}_{\boldsymbol{\mathcal{B}}^{(1)}}]}\left(
\left[
P
\right]_{\varphi(s)}
\right)
\right)
\tag{1}
\\
&=
\mathrm{sc}_{\varphi(s)}^{0[\mathbf{Pth}_{\boldsymbol{\mathcal{B}}^{(1)}}]}\left(
\mathrm{ip}^{([1], Y)@}_{\boldsymbol{\mathcal{B}}^{(1)}, \varphi(s)} \left(
\left[
P
\right]_{\varphi(s)}
\right)
\right)
\tag{2}
\\
&=
\mathrm{sc}_{s}^{0[\mathbf{Pth}_{\boldsymbol{\mathcal{B}}^{(1)}}^{\mathbf{f}^{(1)}}]}\left(
\mathrm{ip}^{([1], Y)@}_{\boldsymbol{\mathcal{B}}^{(1)}, \varphi(s)} \left(
\left[
P
\right]_{\varphi(s)}
\right)
\right).
\tag{3}
\end{align*}
The first equality unravels the interpretation of the operation symbol $\mathrm{sc}_{s}^{0}$ in the partial $\Sigma^{\boldsymbol{\mathcal{A}}^{(1)}}$-algebra $[\mathbf{PT}_{\boldsymbol{\mathcal{B}}^{(1)}}^{\mathbf{f}^{(1)}}]$, introduced in Proposition~\ref{PPTBCatAlg};
the second equality follows from the fact that, according to Claim~\ref{CIsoIpfc}, $\mathrm{ip}^{([1], Y)@}_{\boldsymbol{\mathcal{B}}^{(1)}}$ is a $\Lambda^{\boldsymbol{\mathcal{B}}^{(1)}}$-homomorphism from $[\mathbf{PT}_{\boldsymbol{\mathcal{B}}^{(1)}}]$ to $[\mathbf{Pth}_{\boldsymbol{\mathcal{B}}^{(1)}}]$;
finally, the last equality recovers the interpretation of the operation symbol $\mathrm{sc}_{s}^{0}$ in the partial $\Sigma^{\boldsymbol{\mathcal{A}}^{(1)}}$-algebra $[\mathbf{Pth}_{\boldsymbol{\mathcal{B}}^{(1)}}^{\mathbf{f}^{(1)}}]$, introduced in Proposition~\ref{PQPthBCatAlg}.

Hence, $\mathrm{ip}^{([1], Y)@}_{\boldsymbol{\mathcal{B}}^{(1)}, \varphi}$ is compatible with the $0$-source operation.

{\sffamily The mapping $\mathrm{ip}^{([1], Y)@}_{\boldsymbol{\mathcal{B}}^{(1)}, \varphi}$ is compatible with the $0$-target.}

Let $s$ be a sort in $S$ and let us consider the $0$-target operation symbol $\mathrm{tg}_{s}^{0}$ in $\Sigma^{\boldsymbol{\mathcal{A}}^{(1)}}_{s,s}$. Let $P$ be a path term in $\mathrm{PT}_{\boldsymbol{\mathcal{B}}^{(1)}, \varphi(s)}$, then the following equality holds
$$
\mathrm{ip}^{([1], Y)@}_{\boldsymbol{\mathcal{B}}^{(1)}, \varphi(s)} \left(
\mathrm{tg}_{s}^{0[\mathbf{PT}_{\boldsymbol{\mathcal{B}}^{(1)}}^{\mathbf{f}^{(1)}}]}\left(
\left[
P
\right]_{\varphi(s)}
\right)
\right)
=
\mathrm{tg}_{s}^{0[\mathbf{Pth}_{\boldsymbol{\mathcal{B}}^{(1)}}^{\mathbf{f}^{(1)}}]}\left(
\mathrm{ip}^{([1], Y)@}_{\boldsymbol{\mathcal{B}}^{(1)}, \varphi(s)} \left(
\left[
P
\right]_{\varphi(s)}
\right)
\right).
$$

The proof of this case is identical to that of the $0$-source.

Hence, $\mathrm{ip}^{([1], Y)@}_{\boldsymbol{\mathcal{B}}^{(1)}, \varphi}$ is compatible with the $0$-target operation.

{\sffamily The mapping $\mathrm{ip}^{([1], Y)@}_{\boldsymbol{\mathcal{B}}^{(1)}, \varphi}$ is compatible with the $0$-composition.}

Let $s$ be a sort in $S$ and let us consider the $0$-composition operation symbol $\circ_{s}^{0}$ in $\Sigma_{ss,s}^{\boldsymbol{\mathcal{A}}^{(1)}}$. Let $P$ and $Q$ be two path terms in $\mathrm{PT}_{\boldsymbol{\mathcal{B}}^{(1)}, \varphi(s)}$ such that
$$
\mathrm{sc}_{\boldsymbol{\mathcal{B}}^{(1)}, \varphi(s)}^{(0,1)}\left(
\mathrm{ip}^{(1,Y)@}_{\boldsymbol{\mathcal{B}}^{(1)}, \varphi(s)}\left(
Q
\right)
\right)
=
\mathrm{tg}_{\boldsymbol{\mathcal{B}}^{(1)}, \varphi(s)}^{(0,1)}\left(
\mathrm{ip}^{(1,Y)@}_{\boldsymbol{\mathcal{B}}^{(1)}, \varphi(s)}\left(
P
\right)
\right).
$$
Then the following equality holds
\begin{multline*}
\mathrm{ip}^{([1], Y)@}_{\boldsymbol{\mathcal{B}}^{(1)}, \varphi(s)} \left(
\left[
Q
\right]_{\varphi(s)}
\circ_{s}^{0[\mathbf{PT}_{\boldsymbol{\mathcal{B}}^{(1)}}^{\mathbf{f}^{(1)}}]}
\left[
P
\right]_{\varphi(s)}
\right)
\\
=
\mathrm{ip}^{([1], Y)@}_{\boldsymbol{\mathcal{B}}^{(1)}, \varphi(s)} \left(
\left[
Q
\right]_{\varphi(s)}
\right)
\circ_{s}^{0[\mathbf{Pth}_{\boldsymbol{\mathcal{B}}^{(1)}}^{\mathbf{f}^{(1)}}]}
\mathrm{ip}^{([1], Y)@}_{\boldsymbol{\mathcal{B}}^{(1)}, \varphi(s)} \left(
\left[
P
\right]_{\varphi(s)}
\right).
\end{multline*}

The proof of this case is identical to that of the $0$-source.

Hence, $\mathrm{ip}^{([1], Y)@}_{\boldsymbol{\mathcal{B}}^{(1)}, \varphi}$ is compatible with the $0$-composition operation.

This completes the proof.
\end{proof}

\begin{theorem}\label{TPthBPTB}
The partial $\Sigma^{\boldsymbol{\mathcal{A}}^{(1)}}$-algebras $[\mathbf{Pth}_{\boldsymbol{\mathcal{B}}^{(1)}}^{\mathbf{f}^{(1)}}]$ and $[\mathbf{PT}_{\boldsymbol{\mathcal{B}}^{(1)}}^{\mathbf{f}^{(1)}}]$ are isomorphic.
\end{theorem}

\begin{proof}
Let $s$ be a sort in $S$ and $\mathfrak{P}$ be a path in $\mathrm{Pth}_{\boldsymbol{\mathcal{B}}^{(1)}, \varphi(s)}$.

The following chain of equalities holds
\allowdisplaybreaks
\begin{align*}
\mathrm{ip}^{([1], Y)@}_{\boldsymbol{\mathcal{B}}^{(1)}, \varphi(s)} \left(
\mathrm{CH}^{[1]}_{\boldsymbol{\mathcal{B}}^{(1)}, \varphi(s)} \left(
\left[
\mathfrak{P}
\right]_{\varphi(s)}
\right)
\right)
&=
\left[
\mathrm{ip}^{(2, Y)@}_{\boldsymbol{\mathcal{B}}^{(1)}, \varphi(s)}\left(
\mathrm{CH}^{(1)}_{\boldsymbol{\mathcal{B}}^{(1)}, \varphi(s)}\left(
\mathfrak{P}
\right)
\right)
\right]_{\varphi(s)}
\tag{1}
\\
&=
\left[
\mathfrak{P}
\right]_{\varphi(s)}.
\tag{2}
\end{align*}
The first equality unravels the definition of the mappings $\mathrm{ip}^{([1], Y)@}_{\boldsymbol{\mathcal{B}}^{(1)}, \varphi}$ and $\mathrm{CH}^{[1]}_{\boldsymbol{\mathcal{B}}^{(1)}, \varphi}$ according to, respectively, Definitions~\ref{DPTQIp} and \ref{DPTQCH};
finally, the second equality follows from Proposition~\ref{PIpCH}.

On the other hand, let $s$ be a sort in $S$ and $P$ a path term in $\mathrm{PT}_{\boldsymbol{\mathcal{B}}^{(1)}, \varphi(s)}$.

The following chain of equalities holds
\allowdisplaybreaks
\begin{align*}
\mathrm{CH}^{[1]}_{\boldsymbol{\mathcal{B}}^{(1)}, \varphi(s)} \left(
\mathrm{ip}^{([1], Y)@}_{\boldsymbol{\mathcal{B}}^{(1)}, \varphi(s)} \left(
\left[
P
\right]_{\varphi(s)}
\right)
\right)
&=
\left[
\mathrm{CH}^{(1)}_{\boldsymbol{\mathcal{B}}^{(1)}, \varphi(s)}\left(
\mathrm{ip}^{(2, Y)@}_{\boldsymbol{\mathcal{B}}^{(1)}, \varphi(s)}\left(
P
\right)
\right)
\right]_{\varphi(s)}
\tag{1}
\\
&=
\left[
P
\right]_{\varphi(s)}.
\tag{2}
\end{align*}
The first equality unravels the definition of the mappings $\mathrm{ip}^{([1], Y)@}_{\boldsymbol{\mathcal{B}}^{(1)}, \varphi}$ and $\mathrm{CH}^{[1]}_{\boldsymbol{\mathcal{B}}^{(1)}, \varphi}$ according to, respectively, Definitions~\ref{DPTQIp} and \ref{DPTQCH};
finally, the second equality follows from Proposition~\ref{PIpCH}.
\end{proof}

\section{A structure of partial $\Sigma^{\boldsymbol{\mathcal{A}}^{(1)}}$-algebra on $\T_{\boldsymbol{\mathcal{E}}^{\boldsymbol{\mathcal{B}}^{(1)}}}(\mathbf{Pth}_{\boldsymbol{\mathcal{B}}^{(1)}})$}

We show that $\T_{\boldsymbol{\mathcal{E}}^{\boldsymbol{\mathcal{B}}^{(1)}}}(\mathbf{Pth}_{\boldsymbol{\mathcal{B}}^{(1)}})_{\varphi}$ is equipped, in a natural way with a structure of partial $\Sigma^{\boldsymbol{\mathcal{A}}^{(1)}}$-algebra.

\begin{proposition}
\label{PFreeBCatAlg}
Let $\mathbf{f}^{(1)}=(\varphi, c, (f^{(i)})_{i\in 2})$ be a first-order morphism from $\boldsymbol{\mathcal{A}}^{(1)}$ to $\boldsymbol{\mathcal{B}}^{(1)}$. Then the $S$-sorted set $\T_{\boldsymbol{\mathcal{E}}^{\boldsymbol{\mathcal{B}}^{(1)}}}(\mathbf{Pth}_{\boldsymbol{\mathcal{B}}^{(1)}})_{\varphi}$ is equipped, in a natural way, with a structure of partial $\Sigma^{\boldsymbol{\mathcal{A}}^{(1)}}$-algebra.
\end{proposition}

\begin{proof}
Let us denote by $\mathbf{T}_{\boldsymbol{\mathcal{E}}^{\boldsymbol{\mathcal{B}}^{(1)}}}^{\mathbf{f}^{(1)}}(\mathbf{Pth}_{\boldsymbol{\mathcal{B}}^{(1)}})$ the $\Sigma^{\boldsymbol{\mathcal{A}}^{(1)}}$-algebra defined as follows:

\textsf{(1)}
The underlying $S$-sorted set of $\mathbf{T}_{\boldsymbol{\mathcal{E}}^{\boldsymbol{\mathcal{B}}^{(1)}}}^{\mathbf{f}^{(1)}}(\mathbf{Pth}_{\boldsymbol{\mathcal{B}}^{(1)}})$ is $\T_{\boldsymbol{\mathcal{E}}^{\boldsymbol{\mathcal{B}}^{(1)}}}(\mathbf{Pth}_{\boldsymbol{\mathcal{B}}^{(1)}})_{\varphi}=(\T_{\boldsymbol{\mathcal{E}}^{\boldsymbol{\mathcal{B}}^{(1)}}}(\mathbf{Pth}_{\boldsymbol{\mathcal{B}}^{(1)}})_{\varphi(s)})_{s\in S}$.

\textsf{(2)}
For every $(\mathbf{s}, s)\in S^{\ast}\times S$ and every operation symbol $\sigma\in\Sigma_{\mathbf{s},s}$, the operation $\sigma^{\mathbf{T}_{\boldsymbol{\mathcal{E}}^{\boldsymbol{\mathcal{B}}^{(1)}}}^{\mathbf{f}^{(1)}}(\mathbf{Pth}_{\boldsymbol{\mathcal{B}}^{(1)}})}$ is given by the interpretation of $\sigma$ in the $\Sigma$-algebra $\mathbf{c}_{\mathfrak{d}}^{\ast}(\mathbf{T}_{\boldsymbol{\mathcal{E}}^{\boldsymbol{\mathcal{B}}^{(1)}}}^{(0,1)}(\mathbf{Pth}_{\boldsymbol{\mathcal{B}}^{(1)}}))$. That is, its interpretation is given by the derived operation in $\mathbf{T}_{\boldsymbol{\mathcal{E}}^{\boldsymbol{\mathcal{B}}^{(1)}}}^{(0,1)}(\mathbf{Pth}_{\boldsymbol{\mathcal{B}}^{(1)}})$, the $\Lambda$-reduct of the $\Lambda^{\boldsymbol{\mathcal{B}}^{(1)}}$-algebra  $\mathbf{T}_{\boldsymbol{\mathcal{E}}^{\boldsymbol{\mathcal{B}}^{(1)}}}(\mathbf{Pth}_{\boldsymbol{\mathcal{B}}^{(1)}})$, introduced in Definition~\ref{DVarAbv}.

\textsf{(3)}
For every $s\in S$ and every $\mathfrak{p}\in\mathcal{A}_{s}^{(1)}$, the constant $\mathfrak{p}^{\mathbf{T}_{\boldsymbol{\mathcal{E}}^{\boldsymbol{\mathcal{B}}^{(1)}}}^{\mathbf{f}^{(1)}}(\mathbf{Pth}_{\boldsymbol{\mathcal{B}}^{(1)}})}$ is given by
$$
\mathfrak{p}^{\mathbf{T}_{\boldsymbol{\mathcal{E}}^{\boldsymbol{\mathcal{B}}^{(1)}}}^{\mathbf{f}^{(1)}}(\mathbf{Pth}_{\boldsymbol{\mathcal{B}}^{(1)}})}
=
\mathrm{pr}^{\equiv^{[1]}}_{\boldsymbol{\mathcal{B}}^{(1)}, \varphi(s)} \circ \mathrm{ip}^{(1,Y)@}_{\boldsymbol{\mathcal{B}}^{(1)}, \varphi(s)} \circ \mathrm{CH}^{(1)}_{\boldsymbol{\mathcal{B}}^{(1)}, \varphi(s)} \left(
f^{(1)\flat}_{s}\left(
\mathfrak{p}^{\mathbf{Pth}_{\boldsymbol{\mathcal{A}}^{(1)}}}
\right)
\right).
$$

\textsf{(4)}
For every $s\in S$, the interpretations of the operations $\mathrm{sc}_{s}^{0}$ and $\mathrm{tg}_{s}^{0}$ are given by
\begin{align*}
\mathrm{sc}_{s}^{0\mathbf{T}_{\boldsymbol{\mathcal{E}}^{\boldsymbol{\mathcal{B}}^{(1)}}}^{\mathbf{f}^{(1)}}(\mathbf{Pth}_{\boldsymbol{\mathcal{B}}^{(1)}})}
&=
\mathrm{sc}_{\varphi(s)}^{0\mathbf{T}_{\boldsymbol{\mathcal{E}}^{\boldsymbol{\mathcal{B}}^{(1)}}}(\mathbf{Pth}_{\boldsymbol{\mathcal{B}}^{(1)}})}
&&\mbox{and}&
\mathrm{tg}_{s}^{0\mathbf{T}_{\boldsymbol{\mathcal{E}}^{\boldsymbol{\mathcal{B}}^{(1)}}}^{\mathbf{f}^{(1)}}(\mathbf{Pth}_{\boldsymbol{\mathcal{B}}^{(1)}})}
&=
\mathrm{tg}_{\varphi(s)}^{0\mathbf{T}_{\boldsymbol{\mathcal{E}}^{\boldsymbol{\mathcal{B}}^{(1)}}}(\mathbf{Pth}_{\boldsymbol{\mathcal{B}}^{(1)}})}.
\end{align*}
That is, their interpretations are given by the interpretations of $\mathrm{sc}_{\varphi(s)}^{0}$ and $\mathrm{tg}_{\varphi(s)}^{0}$ in $\mathbf{T}_{\boldsymbol{\mathcal{E}}^{\boldsymbol{\mathcal{B}}^{(1)}}}(\mathbf{Pth}_{\boldsymbol{\mathcal{B}}^{(1)}})$ that, we recall, was introduced in Definition~\ref{DVarAbv}.

\textsf{(5)}
Similarly,  for every $s \in S$, the interpretations of the operation $\circ_{s}^{0}$ is given by
$$
\circ_{s}^{0\mathbf{T}_{\boldsymbol{\mathcal{E}}^{\boldsymbol{\mathcal{B}}^{(1)}}}^{\mathbf{f}^{(1)}}(\mathbf{Pth}_{\boldsymbol{\mathcal{B}}^{(1)}})}
=
\circ_{\varphi(s)}^{0\mathbf{T}_{\boldsymbol{\mathcal{E}}^{\boldsymbol{\mathcal{B}}^{(1)}}}(\mathbf{Pth}_{\boldsymbol{\mathcal{B}}^{(1)}})}.
$$
That is, its interpretation is given by the interpretation of $\circ_{\varphi(s)}^{0}$ in $\mathbf{T}_{\boldsymbol{\mathcal{E}}^{\boldsymbol{\mathcal{B}}^{(1)}}}(\mathbf{Pth}_{\boldsymbol{\mathcal{B}}^{(1)}})$ that, we recall, was introduced in Definition~\ref{DVarAbv}.

This completes the definition of the partial $\Sigma^{\boldsymbol{\mathcal{A}}^{(1)}}$-algebra $\mathbf{T}_{\boldsymbol{\mathcal{E}}^{\boldsymbol{\mathcal{B}}^{(1)}}}^{\mathbf{f}^{(1)}}(\mathbf{Pth}_{\boldsymbol{\mathcal{B}}^{(1)}})$.
\end{proof}

Finally, we show that the partial $\Sigma^{\boldsymbol{\mathcal{A}}^{(1)}}$-algebras $[\mathbf{Pth}_{\boldsymbol{\mathcal{B}}^{(1)}}^{\mathbf{f}^{(1)}}]_{\varphi}$ and $\mathbf{T}_{\boldsymbol{\mathcal{E}}^{\boldsymbol{\mathcal{B}}^{(1)}}}^{\mathbf{f}^{(1)}}(\mathbf{Pth}_{\boldsymbol{\mathcal{B}}^{(1)}})$ are isomorphic. In virtue of Theorem~\ref{TPthBPTB} the same relation will apply to the partial $\Sigma^{\boldsymbol{\mathcal{A}}^{(1)}}$-algebra $[ \mathbf{PT}^{\mathbf{f}^{(1)}}_{\boldsymbol{\mathcal{B}}^{(1)}} ]$.

\begin{proposition}
\label{PQVarCatHom}
Let $\mathbf{f}^{(1)}=(\varphi, c, (f^{(i)})_{i\in 2})$ be a first-order morphism from $\boldsymbol{\mathcal{A}}^{(1)}$ to $\boldsymbol{\mathcal{B}}^{(1)}$. Then the mapping 
$$
(\mathrm{pr}^{\equiv^{[1]}}_{\boldsymbol{\mathcal{B}}^{(1)}} \circ \mathrm{ip}^{(1,Y)@}_{\boldsymbol{\mathcal{B}}^{(1)}} \circ \mathrm{CH}^{(1)}_{\boldsymbol{\mathcal{B}}^{(1)}})^{\natural}_{\varphi}
\colon 
[\mathrm{Pth}_{\boldsymbol{\mathcal{B}}^{(1)}}]_{\varphi}
\mor
\T_{\boldsymbol{\mathcal{E}}^{\boldsymbol{\mathcal{B}}^{(1)}}}(\mathbf{Pth}_{\boldsymbol{\mathcal{B}}^{(1)}})_{\varphi}
$$
is a $\Sigma^{\boldsymbol{\mathcal{A}}^{(1)}}$-homomorphism from $[\mathbf{Pth}_{\boldsymbol{\mathcal{B}}^{(1)}}^{\mathbf{f}^{(1)}}]$ to $\mathbf{T}_{\boldsymbol{\mathcal{E}}^{\boldsymbol{\mathcal{B}}^{(1)}}}^{\mathbf{f}^{(1)}}(\mathbf{Pth}_{\boldsymbol{\mathcal{B}}^{(1)}})$.
\end{proposition}

\begin{proof}
We prove that $(\mathrm{pr}^{\equiv^{[1]}}_{\boldsymbol{\mathcal{B}}^{(1)}} \circ \mathrm{ip}^{(1,Y)@}_{\boldsymbol{\mathcal{B}}^{(1)}} \circ \mathrm{CH}^{(1)}_{\boldsymbol{\mathcal{B}}^{(1)}})^{\natural}_{\varphi}$ is compatible with every operation symbol in $\Sigma^{\boldsymbol{\mathcal{A}}^{(1)}}$.

{\sffamily The mapping $(\mathrm{pr}^{\equiv^{[1]}}_{\boldsymbol{\mathcal{B}}^{(1)}} \circ \mathrm{ip}^{(1,Y)@}_{\boldsymbol{\mathcal{B}}^{(1)}} \circ \mathrm{CH}^{(1)}_{\boldsymbol{\mathcal{B}}^{(1)}})^{\natural}_{\varphi}$ is a $\Sigma$-homomorphism.}

Note that $(\mathrm{pr}^{\equiv^{[1]}}_{\boldsymbol{\mathcal{B}}^{(1)}} \circ \mathrm{ip}^{(1,Y)@}_{\boldsymbol{\mathcal{B}}^{(1)}} \circ \mathrm{CH}^{(1)}_{\boldsymbol{\mathcal{B}}^{(1)}})^{\natural}_{\varphi} = \mathbf{c}_{\mathfrak{d}}^{\ast} ((\mathrm{pr}^{\equiv^{[1]}}_{\boldsymbol{\mathcal{B}}^{(1)}} \circ \mathrm{ip}^{(1,Y)@}_{\boldsymbol{\mathcal{B}}^{(1)}} \circ \mathrm{CH}^{(1)}_{\boldsymbol{\mathcal{B}}^{(1)}})^{\natural})$. By Remark~\ref{RQPUniv}, the mapping $(\mathrm{pr}^{\equiv^{[1]}}_{\boldsymbol{\mathcal{B}}^{(1)}} \circ \mathrm{ip}^{(1,Y)@}_{\boldsymbol{\mathcal{B}}^{(1)}} \circ \mathrm{CH}^{(1)}_{\boldsymbol{\mathcal{B}}^{(1)}})^{\natural}$ is a $\Lambda^{\boldsymbol{\mathcal{B}}^{(1)}}$-homomorphism, thus in particular a $\Lambda$-homomorphism. Therefore, it follows from Proposition~\ref{PFunSig} that the mapping $(\mathrm{pr}^{\equiv^{[1]}}_{\boldsymbol{\mathcal{B}}^{(1)}} \circ \mathrm{ip}^{(1,Y)@}_{\boldsymbol{\mathcal{B}}^{(1)}} \circ \mathrm{CH}^{(1)}_{\boldsymbol{\mathcal{B}}^{(1)}})^{\natural}_{\varphi}$ is a $\Sigma$-homomorphism.

{\sffamily The mapping $(\mathrm{pr}^{\equiv^{[1]}}_{\boldsymbol{\mathcal{B}}^{(1)}} \circ \mathrm{ip}^{(1,Y)@}_{\boldsymbol{\mathcal{B}}^{(1)}} \circ \mathrm{CH}^{(1)}_{\boldsymbol{\mathcal{B}}^{(1)}})^{\natural}_{\varphi}$ is compatible with the rewrite rules.}

Let $s$ be a sort in $S$ and $\mathfrak{p}$ a rewrite rule in $\mathcal{A}_{s}^{(1)}$. Thus, the following chain of equalities holds
\begin{flushleft}
$
(\mathrm{pr}^{\equiv^{[1]}}_{\boldsymbol{\mathcal{B}}^{(1)}} \circ \mathrm{ip}^{(1,Y)@}_{\boldsymbol{\mathcal{B}}^{(1)}} \circ \mathrm{CH}^{(1)}_{\boldsymbol{\mathcal{B}}^{(1)}})^{\natural}_{\varphi(s)} \left(
\mathfrak{p}^{[\mathbf{Pth}_{\boldsymbol{\mathcal{B}}^{(1)}}^{\mathbf{f}^{(1)}}]}
\right)
$
\allowdisplaybreaks
\begin{align*}
&=
(\mathrm{pr}^{\equiv^{[1]}}_{\boldsymbol{\mathcal{B}}^{(1)}} \circ \mathrm{ip}^{(1,Y)@}_{\boldsymbol{\mathcal{B}}^{(1)}} \circ \mathrm{CH}^{(1)}_{\boldsymbol{\mathcal{B}}^{(1)}})^{\natural}_{\varphi(s)} \left(
\left[
f^{(1)\flat}\left(
\mathfrak{p}^{\mathbf{Pth}_{\boldsymbol{\mathcal{A}}^{(1)}}}
\right)
\right]_{\varphi(s)}
\right)
\tag{1}
\\
&=
(\mathrm{pr}^{\equiv^{[1]}}_{\boldsymbol{\mathcal{B}}^{(1)}} \circ \mathrm{ip}^{(1,Y)@}_{\boldsymbol{\mathcal{B}}^{(1)}} \circ \mathrm{CH}^{(1)}_{\boldsymbol{\mathcal{B}}^{(1)}})_{\varphi(s)} \left(
f^{(1)\flat}\left(
\mathfrak{p}^{\mathbf{Pth}_{\boldsymbol{\mathcal{A}}^{(1)}}}
\right)
\right)
\tag{2}
\\
&=
\mathrm{pr}^{\equiv^{[1]}}_{\boldsymbol{\mathcal{B}}^{(1)}, \varphi(s)} 
\circ
\mathrm{ip}^{(1,Y)@}_{\boldsymbol{\mathcal{B}}^{(1)}, \varphi(s)}
\circ
\mathrm{CH}^{(1)}_{\boldsymbol{\mathcal{B}}^{(1)}, \varphi(s)} \left(
f^{(1)\flat}\left(
\mathfrak{p}^{\mathbf{Pth}_{\boldsymbol{\mathcal{A}}^{(1)}}}
\right)
\right)
\tag{3}
\\
&=
\mathfrak{p}^{\mathbf{T}_{\boldsymbol{\mathcal{E}}^{\boldsymbol{\mathcal{B}}^{(1)}}}^{\mathbf{f}^{(1)}}(\mathbf{Pth}_{\boldsymbol{\mathcal{B}}^{(1)}})}.
\tag{4}
\end{align*}
\end{flushleft}

The first equality unravels the interpretation of the constant $\mathfrak{p}$ in the partial $\Sigma^{\boldsymbol{\mathcal{A}}^{(1)}}$-algebra $[\mathbf{Pth}^{\mathbf{f}^{(1)}}_{\boldsymbol{\mathcal{B}}^{(1)}}]$, introduced in Proposition~\ref{PQPthBCatAlg};
the second equality follows from Remark~\ref{RQPUniv};
the third equality unravels the definition of the $\varphi(s)$ component of a composition of many-sorted mappings;
finally, the last equality recovers the interpretation of the constant $\mathfrak{p}$ in the partial $\Sigma^{\boldsymbol{\mathcal{A}}^{(1)}}$-algebra $\mathbf{T}_{\boldsymbol{\mathcal{E}}^{\boldsymbol{\mathcal{B}}^{(1)}}}^{\mathbf{f}^{(1)}}(\mathbf{Pth}_{\boldsymbol{\mathcal{B}}^{(1)}})$, introduced in Proposition~\ref{PFreeBCatAlg}.

Hence, $(\mathrm{pr}^{\equiv^{[1]}}_{\boldsymbol{\mathcal{B}}^{(1)}} \circ \mathrm{ip}^{(1,Y)@}_{\boldsymbol{\mathcal{B}}^{(1)}} \circ \mathrm{CH}^{(1)}_{\boldsymbol{\mathcal{B}}^{(1)}})^{\natural}_{\varphi}$ is compatible with the rewrite rules.

{\sffamily The mapping $(\mathrm{pr}^{\equiv^{[1]}}_{\boldsymbol{\mathcal{B}}^{(1)}} \circ \mathrm{ip}^{(1,Y)@}_{\boldsymbol{\mathcal{B}}^{(1)}} \circ \mathrm{CH}^{(1)}_{\boldsymbol{\mathcal{B}}^{(1)}})^{\natural}_{\varphi}$ is compatible with the $0$-source.}

Let $s$ be a sort in $S$ and let us consider the $0$-source operation symbol $\mathrm{sc}_{s}^{0}$ in $\Sigma^{\boldsymbol{\mathcal{A}}^{(1)}}_{s,s}$. Let $\mathfrak{P}$ be a path in $\mathrm{Pth}_{\boldsymbol{\mathcal{B}}^{(1)}, \varphi(s)}$.

The following chain of equalities holds
\begin{flushleft}
$
(\mathrm{pr}^{\equiv^{[1]}}_{\boldsymbol{\mathcal{B}}^{(1)}} \circ \mathrm{ip}^{(1,Y)@}_{\boldsymbol{\mathcal{B}}^{(1)}} \circ \mathrm{CH}^{(1)}_{\boldsymbol{\mathcal{B}}^{(1)}})^{\natural}_{\varphi(s)} \left(
\mathrm{sc}_{s}^{0[\mathbf{Pth}_{\boldsymbol{\mathcal{B}}^{(1)}}^{\mathbf{f}^{(1)}}]}\left(
\left[
\mathfrak{P}
\right]_{\varphi(s)}
\right)
\right)
$
\allowdisplaybreaks
\begin{align*}
&=
(\mathrm{pr}^{\equiv^{[1]}}_{\boldsymbol{\mathcal{B}}^{(1)}} \circ \mathrm{ip}^{(1,Y)@}_{\boldsymbol{\mathcal{B}}^{(1)}} \circ \mathrm{CH}^{(1)}_{\boldsymbol{\mathcal{B}}^{(1)}})^{\natural}_{\varphi(s)} \left(
\mathrm{sc}_{\varphi(s)}^{0[\mathbf{Pth}_{\boldsymbol{\mathcal{B}}^{(1)}}]}\left(
\left[
\mathfrak{P}
\right]_{\varphi(s)}
\right)
\right)
\tag{1}
\\
&=
(\mathrm{pr}^{\equiv^{[1]}}_{\boldsymbol{\mathcal{B}}^{(1)}} \circ \mathrm{ip}^{(1,Y)@}_{\boldsymbol{\mathcal{B}}^{(1)}} \circ \mathrm{CH}^{(1)}_{\boldsymbol{\mathcal{B}}^{(1)}})^{\natural}_{\varphi(s)} \left(
\left[
\mathrm{sc}_{\varphi(s)}^{0\mathbf{Pth}_{\boldsymbol{\mathcal{B}}^{(1)}}}\left(
\mathfrak{P}
\right)
\right]_{\varphi(s)}
\right)
\tag{2}
\\
&=
(\mathrm{pr}^{\equiv^{[1]}}_{\boldsymbol{\mathcal{B}}^{(1)}} \circ \mathrm{ip}^{(1,Y)@}_{\boldsymbol{\mathcal{B}}^{(1)}} \circ \mathrm{CH}^{(1)}_{\boldsymbol{\mathcal{B}}^{(1)}})_{\varphi(s)} \left(
\mathrm{sc}_{\varphi(s)}^{0\mathbf{Pth}_{\boldsymbol{\mathcal{B}}^{(1)}}}\left(
\mathfrak{P}
\right)
\right)
\tag{3}
\\
&=
\mathrm{sc}_{\varphi(s)}^{0\mathbf{T}_{\boldsymbol{\mathcal{E}}^{\boldsymbol{\mathcal{B}}^{(1)}}}^{\mathbf{f}^{(1)}}(\mathbf{Pth}_{\boldsymbol{\mathcal{B}}^{(1)}})}\left(
(\mathrm{pr}^{\equiv^{[1]}}_{\boldsymbol{\mathcal{B}}^{(1)}} \circ \mathrm{ip}^{(1,Y)@}_{\boldsymbol{\mathcal{B}}^{(1)}} \circ \mathrm{CH}^{(1)}_{\boldsymbol{\mathcal{B}}^{(1)}})_{\varphi(s)} \left(
\mathfrak{P}
\right)
\right)
\tag{4}
\\
&=
\mathrm{sc}_{\varphi(s)}^{0\mathbf{T}_{\boldsymbol{\mathcal{E}}^{\boldsymbol{\mathcal{B}}^{(1)}}}^{\mathbf{f}^{(1)}}(\mathbf{Pth}_{\boldsymbol{\mathcal{B}}^{(1)}})}\left(
(\mathrm{pr}^{\equiv^{[1]}}_{\boldsymbol{\mathcal{B}}^{(1)}} \circ \mathrm{ip}^{(1,Y)@}_{\boldsymbol{\mathcal{B}}^{(1)}} \circ \mathrm{CH}^{(1)}_{\boldsymbol{\mathcal{B}}^{(1)}})^{\natural}_{\varphi(s)} \left(
\left[
\mathfrak{P}
\right]_{\varphi(s)}
\right)
\right)
\tag{5}
\\
&=
\mathrm{sc}_{s}^{0\mathbf{T}_{\boldsymbol{\mathcal{E}}^{\boldsymbol{\mathcal{B}}^{(1)}}}^{\mathbf{f}^{(1)}}(\mathbf{Pth}_{\boldsymbol{\mathcal{B}}^{(1)}})}\left(
(\mathrm{pr}^{\equiv^{[1]}}_{\boldsymbol{\mathcal{B}}^{(1)}} \circ \mathrm{ip}^{(1,Y)@}_{\boldsymbol{\mathcal{B}}^{(1)}} \circ \mathrm{CH}^{(1)}_{\boldsymbol{\mathcal{B}}^{(1)}})^{\natural}_{\varphi(s)} \left(
\left[
\mathfrak{P}
\right]_{\varphi(s)}
\right)
\right)
\tag{6}
\end{align*}
\end{flushleft}

The first equality unravels the interpretation of the operation symbol $\mathrm{sc}_{s}^{0}$ in the partial $\Sigma^{\boldsymbol{\mathcal{A}}^{(1)}}$-algebra $[\mathbf{Pth}_{\boldsymbol{\mathcal{B}}^{(1)}}^{\mathbf{f}^{(1)}}]$, introduced in Proposition~\ref{PQPthBCatAlg};
the second equality unravels the interpretation of the operation symbol $\mathrm{sc}_{s}^{0}$ in the partial $\Sigma^{\boldsymbol{\mathcal{A}}^{(1)}}$-algebra $[\mathbf{Pth}_{\boldsymbol{\mathcal{B}}^{(1)}}]$, introduced in Proposition~\ref{PCHCatAlg};
the third equality follows from Remark~\ref{RQPUniv};
the fourth follows from the fact that, according to Proposition~\ref{PVarKer}, the mapping $(\mathrm{pr}^{\equiv^{[1]}}_{\boldsymbol{\mathcal{B}}^{(1)}} \circ \mathrm{ip}^{(1,Y)@}_{\boldsymbol{\mathcal{B}}^{(1)}} \circ \mathrm{CH}^{(1)}_{\boldsymbol{\mathcal{B}}^{(1)}})_{\varphi(s)}$ is a $\Lambda^{\boldsymbol{\mathcal{B}}^{(1)}}$-homomorphism;
the fifth equality follows from Remark~	\ref{RQPUniv};
finally, the last equality recovers the interpretation of the operation symbol $\mathrm{sc}_{s}^{0}$ in the partial $\Sigma^{\boldsymbol{\mathcal{A}}^{(1)}}$-algebra $\mathbf{T}_{\boldsymbol{\mathcal{E}}^{\boldsymbol{\mathcal{B}}^{(1)}}}^{\mathbf{f}^{(1)}}(\mathbf{Pth}_{\boldsymbol{\mathcal{B}}^{(1)}})$, introduced in Proposition~\ref{PFreeBCatAlg}.

Hence, $(\mathrm{pr}^{\equiv^{[1]}}_{\boldsymbol{\mathcal{B}}^{(1)}} \circ \mathrm{ip}^{(1,Y)@}_{\boldsymbol{\mathcal{B}}^{(1)}} \circ \mathrm{CH}^{(1)}_{\boldsymbol{\mathcal{B}}^{(1)}})^{\natural}_{\varphi}$ is compatible with the $0$-source operation.

{\sffamily The mapping $(\mathrm{pr}^{\equiv^{[1]}}_{\boldsymbol{\mathcal{B}}^{(1)}} \circ \mathrm{ip}^{(1,Y)@}_{\boldsymbol{\mathcal{B}}^{(1)}} \circ \mathrm{CH}^{(1)}_{\boldsymbol{\mathcal{B}}^{(1)}})^{\natural}_{\varphi}$ is compatible with the $0$-target.}

Let $s$ be a sort in $S$ and let us consider the $0$-target operation symbol $\mathrm{tg}_{s}^{0}$ in $\Sigma^{\boldsymbol{\mathcal{A}}^{(1)}}_{s,s}$. Let $\mathfrak{P}$ be a path in $\mathrm{Pth}_{\boldsymbol{\mathcal{B}}^{(1)}, \varphi(s)}$, then the following equality holds
\begin{multline*}
(\mathrm{pr}^{\equiv^{[1]}}_{\boldsymbol{\mathcal{B}}^{(1)}} \circ \mathrm{ip}^{(1,Y)@}_{\boldsymbol{\mathcal{B}}^{(1)}} \circ \mathrm{CH}^{(1)}_{\boldsymbol{\mathcal{B}}^{(1)}})^{\natural}_{\varphi(s)} \left(
\mathrm{tg}_{s}^{0[\mathbf{Pth}_{\boldsymbol{\mathcal{B}}^{(1)}}^{\mathbf{f}^{(1)}}]}\left(
\left[
\mathfrak{P}
\right]_{\varphi(s)}
\right)
\right)
\\
=
\mathrm{tg}_{\varphi(s)}^{0\mathbf{T}_{\boldsymbol{\mathcal{E}}^{\boldsymbol{\mathcal{B}}^{(1)}}}^{\mathbf{f}^{(1)}}(\mathbf{Pth}_{\boldsymbol{\mathcal{B}}^{(1)}})}\left(
(\mathrm{pr}^{\equiv^{[1]}}_{\boldsymbol{\mathcal{B}}^{(1)}} \circ \mathrm{ip}^{(1,Y)@}_{\boldsymbol{\mathcal{B}}^{(1)}} \circ \mathrm{CH}^{(1)}_{\boldsymbol{\mathcal{B}}^{(1)}})^{\natural}_{\varphi(s)} \left(
\left[
\mathfrak{P}
\right]_{\varphi(s)}
\right)
\right).
\end{multline*}

The proof of this case is identical to that of the $0$-source.

Hence, $(\mathrm{pr}^{\equiv^{[1]}}_{\boldsymbol{\mathcal{B}}^{(1)}} \circ \mathrm{ip}^{(1,Y)@}_{\boldsymbol{\mathcal{B}}^{(1)}} \circ \mathrm{CH}^{(1)}_{\boldsymbol{\mathcal{B}}^{(1)}})^{\natural}_{\varphi}$ is compatible with the $0$-target operation.

{\sffamily The mapping $(\mathrm{pr}^{\equiv^{[1]}}_{\boldsymbol{\mathcal{B}}^{(1)}} \circ \mathrm{ip}^{(1,Y)@}_{\boldsymbol{\mathcal{B}}^{(1)}} \circ \mathrm{CH}^{(1)}_{\boldsymbol{\mathcal{B}}^{(1)}})^{\natural}_{\varphi}$ is compatible with the $0$-composition.}

Let $s$ be a sort in $S$ and let us consider the $0$-composition operation symbol $\circ_{s}^{0}$ in $\Sigma_{ss,s}^{\boldsymbol{\mathcal{A}}^{(1)}}$. Let $\mathfrak{P}$ and $\mathfrak{Q}$ be two paths in $\mathrm{Pth}_{\boldsymbol{\mathcal{B}}^{(1)}, \varphi(s)}$ such that
$$
\mathrm{sc}_{\boldsymbol{\mathcal{B}}^{(1)}, \varphi(s)}^{(0,1)}\left(\mathfrak{Q}\right)
=
\mathrm{tg}_{\boldsymbol{\mathcal{B}}^{(1)}, \varphi(s)}^{(0,1)}\left(\mathfrak{P}\right).
$$

Then the following equality holds
\allowdisplaybreaks
\begin{align*}
&
(\mathrm{pr}^{\equiv^{[1]}}_{\boldsymbol{\mathcal{B}}^{(1)}} \circ \mathrm{ip}^{(1,Y)@}_{\boldsymbol{\mathcal{B}}^{(1)}} \circ \mathrm{CH}^{(1)}_{\boldsymbol{\mathcal{B}}^{(1)}})^{\natural}_{\varphi(s)} \left(
\left[
\mathfrak{Q}
\right]_{\varphi(s)}
\circ_{s}^{0[\mathbf{Pth}_{\boldsymbol{\mathcal{B}}^{(1)}}^{\mathbf{f}^{(1)}}]}
\left[
\mathfrak{P}
\right]_{\varphi(s)}
\right)
\\
&\hspace{1.5cm}
=(\mathrm{pr}^{\equiv^{[1]}}_{\boldsymbol{\mathcal{B}}^{(1)}} \circ \mathrm{ip}^{(1,Y)@}_{\boldsymbol{\mathcal{B}}^{(1)}} \circ \mathrm{CH}^{(1)}_{\boldsymbol{\mathcal{B}}^{(1)}})^{\natural}_{\varphi(s)} \left(
\left[
\mathfrak{Q}
\right]_{\varphi(s)}
\right)
\\
&\hspace{2cm}
\circ_{s}^{0\mathbf{T}_{\boldsymbol{\mathcal{E}}^{\boldsymbol{\mathcal{B}}^{(1)}}}^{\mathbf{f}^{(1)}}(\mathbf{Pth}_{\boldsymbol{\mathcal{B}}^{(1)}})}
(\mathrm{pr}^{\equiv^{[1]}}_{\boldsymbol{\mathcal{B}}^{(1)}} \circ \mathrm{ip}^{(1,Y)@}_{\boldsymbol{\mathcal{B}}^{(1)}} \circ \mathrm{CH}^{(1)}_{\boldsymbol{\mathcal{B}}^{(1)}})^{\natural}_{\varphi(s)} \left(
\left[
\mathfrak{P}
\right]_{\varphi(s)}
\right).
\end{align*}

The proof of this case is identical to that of the $0$-source.

Hence, $(\mathrm{pr}^{\equiv^{[1]}}_{\boldsymbol{\mathcal{B}}^{(1)}} \circ \mathrm{ip}^{(1,Y)@}_{\boldsymbol{\mathcal{B}}^{(1)}} \circ \mathrm{CH}^{(1)}_{\boldsymbol{\mathcal{B}}^{(1)}})^{\natural}_{\varphi}$ is compatible with the $0$-composition operation.

This completes the proof.
\end{proof}

\begin{proposition}
\label{PFreeBQCHBCatHom}
Let $\mathbf{f}^{(1)}=(\varphi, c, (f^{(i)})_{i\in 2})$ be a first-order morphism from $\boldsymbol{\mathcal{A}}^{(1)}$ to $\boldsymbol{\mathcal{B}}^{(1)}$. Then the mapping 
$$
\mathrm{pr}^{\mathrm{Ker}(\mathrm{CH}^{(1)})\mathsf{p}}_{\boldsymbol{\mathcal{B}}^{(1)}, \varphi}
\colon 
\T_{\boldsymbol{\mathcal{E}}^{\boldsymbol{\mathcal{B}}^{(1)}}}(\mathbf{Pth}_{\boldsymbol{\mathcal{B}}^{(1)}})_{\varphi}
\mor
[\mathrm{Pth}_{\boldsymbol{\mathcal{B}}^{(1)}}]_{\varphi}
$$
is a $\Sigma^{\boldsymbol{\mathcal{A}}^{(1)}}$-homomorphism from $\mathbf{T}_{\boldsymbol{\mathcal{E}}^{\boldsymbol{\mathcal{B}}^{(1)}}}^{\mathbf{f}^{(1)}}(\mathbf{Pth}_{\boldsymbol{\mathcal{B}}^{(1)}})$ to $[\mathbf{Pth}_{\boldsymbol{\mathcal{B}}^{(1)}}^{\mathbf{f}^{(1)}}]$.
\end{proposition}

\begin{proof}
We prove that $\mathrm{pr}^{\mathrm{Ker}(\mathrm{CH}^{(1)})\mathsf{p}}_{\boldsymbol{\mathcal{B}}^{(1)}, \varphi}$ is compatible with every operation symbol in $\Sigma^{\boldsymbol{\mathcal{A}}^{(1)}}$.

{\sffamily The mapping $\mathrm{pr}^{\mathrm{Ker}(\mathrm{CH}^{(1)})\mathsf{p}}_{\boldsymbol{\mathcal{B}}^{(1)}, \varphi}$ is a $\Sigma$-homomorphism.}

Note that $\mathrm{pr}^{\mathrm{Ker}(\mathrm{CH}^{(1)})\mathsf{p}}_{\boldsymbol{\mathcal{B}}^{(1)}, \varphi} = \mathbf{c}_{\mathfrak{d}}^{\ast} (\mathrm{pr}^{\mathrm{Ker}(\mathrm{CH}^{(1)})\mathsf{p}}_{\boldsymbol{\mathcal{B}}^{(1)}})$. By Corollary~\ref{CVarPr}, the mapping $\mathrm{pr}^{\mathrm{Ker}(\mathrm{CH}^{(1)})\mathsf{p}}_{\boldsymbol{\mathcal{B}}^{(1)}}$ is a $\Lambda^{\boldsymbol{\mathcal{B}}^{(1)}}$-homomorphism, thus in particular a $\Lambda$-homomorphism. Therefore, it follows from Proposition~\ref{PFunSig} that the mapping $\mathrm{pr}^{\mathrm{Ker}(\mathrm{CH}^{(1)})\mathsf{p}}_{\boldsymbol{\mathcal{B}}^{(1)}, \varphi}$ is a $\Sigma$-homomorphism.

{\sffamily The mapping $\mathrm{pr}^{\mathrm{Ker}(\mathrm{CH}^{(1)})\mathsf{p}}_{\boldsymbol{\mathcal{B}}^{(1)}, \varphi}$ is compatible with the rewrite rules.}

Let $s$ be a sort in $S$ and $\mathfrak{p}$ a rewrite rule in $\mathcal{A}_{s}^{(1)}$. Thus, the following chain of equalities holds
\begin{flushleft}
$
\mathrm{pr}^{\mathrm{Ker}(\mathrm{CH}^{(1)})\mathsf{p}}_{\boldsymbol{\mathcal{B}}^{(1)}, \varphi(s)} \left(
\mathfrak{p}^{\mathbf{T}_{\boldsymbol{\mathcal{E}}^{\boldsymbol{\mathcal{B}}^{(1)}}}^{\mathbf{f}^{(1)}}(\mathbf{Pth}_{\boldsymbol{\mathcal{B}}^{(1)}})}
\right)
$
\allowdisplaybreaks
\begin{align*}
&=
\mathrm{pr}^{\mathrm{Ker}(\mathrm{CH}^{(1)})\mathsf{p}}_{\boldsymbol{\mathcal{B}}^{(1)}, \varphi(s)} \left(
\mathrm{pr}^{\equiv^{[1]}}_{\boldsymbol{\mathcal{B}}^{(1)}, \varphi(s)} \circ \mathrm{ip}^{(1,Y)@}_{\boldsymbol{\mathcal{B}}^{(1)}, \varphi(s)} \circ \mathrm{CH}^{(1)}_{\boldsymbol{\mathcal{B}}^{(1)}, \varphi(s)} \left(
f^{(1)\flat}_{s}\left(
\mathfrak{p}^{\mathbf{Pth}_{\boldsymbol{\mathcal{A}}^{(1)}}}
\right)
\right)
\right)
\tag{1}
\\
&=
\mathrm{pr}^{\mathrm{Ker}(\mathrm{CH}^{(1)})\mathsf{p}}_{\boldsymbol{\mathcal{B}}^{(1)}, \varphi(s)} \left(
(\mathrm{pr}^{\equiv^{[1]}}_{\boldsymbol{\mathcal{B}}^{(1)}} \circ \mathrm{ip}^{(1,Y)@}_{\boldsymbol{\mathcal{B}}^{(1)}} \circ \mathrm{CH}^{(1)}_{\boldsymbol{\mathcal{B}}^{(1)}})_{\varphi(s)} \left(
f^{(1)\flat}_{s}\left(
\mathfrak{p}^{\mathbf{Pth}_{\boldsymbol{\mathcal{A}}^{(1)}}}
\right)
\right)
\right)
\tag{2}
\\
&=
\resizebox{.89\textwidth}{!}{%
$
\mathrm{pr}^{\mathrm{Ker}(\mathrm{CH}^{(1)})\mathsf{p}}_{\boldsymbol{\mathcal{B}}^{(1)}, \varphi(s)} \left(
(\mathrm{pr}^{\equiv^{[1]}}_{\boldsymbol{\mathcal{B}}^{(1)}} \circ \mathrm{ip}^{(1,Y)@}_{\boldsymbol{\mathcal{B}}^{(1)}} \circ \mathrm{CH}^{(1)}_{\boldsymbol{\mathcal{B}}^{(1)}})_{\varphi(s)}^{\natural} \left(
\left[
f^{(1)\flat}_{s}\left(
\mathfrak{p}^{\mathbf{Pth}_{\boldsymbol{\mathcal{A}}^{(1)}}}
\right)
\right]_{\varphi(s)}
\right)
\right)
$
}
\tag{3}
\\
&=
\left[
f^{(1)\flat}_{s}\left(
\mathfrak{p}^{\mathbf{Pth}_{\boldsymbol{\mathcal{A}}^{(1)}}}
\right)
\right]_{\varphi(s)}
\tag{4}
\\
&=
\mathfrak{p}^{[\mathbf{Pth}_{\boldsymbol{\mathcal{B}}^{(1)}}^{\mathbf{f}^{(1)}}]}.
\end{align*}
\end{flushleft}

The first equality unravels the interpretation of the constant $\mathfrak{p}$ in the partial $\Sigma^{\boldsymbol{\mathcal{A}}^{(1)}}$-algebra $\mathbf{T}_{\boldsymbol{\mathcal{E}}^{\boldsymbol{\mathcal{B}}^{(1)}}}^{\mathbf{f}^{(1)}}(\mathbf{Pth}_{\boldsymbol{\mathcal{B}}^{(1)}})$, introduced in Proposition~\ref{PFreeBCatAlg};
the second equality recovers the definition of the $\varphi(s)$ component of a composition of many-sorted mappings;
the third equality follows from Remark~\ref{RQPUniv};
The fourth equality follows from Theorem~\ref{TPthFree};
finally, the last equality recovers the interpretation of the constant $\mathfrak{p}$ in the partial $\Sigma^{\boldsymbol{\mathcal{A}}^{(1)}}$-algebra $[\mathbf{Pth}_{\boldsymbol{\mathcal{B}}^{(1)}}^{\mathbf{f}^{(1)}}]$, introduced in Proposition~\ref{PQPthBCatAlg}.

Hence, $\mathrm{pr}^{\mathrm{Ker}(\mathrm{CH}^{(1)})\mathsf{p}}_{\boldsymbol{\mathcal{B}}^{(1)}, \varphi}$ is compatible with the rewrite rules.

{\sffamily The mapping $\mathrm{pr}^{\mathrm{Ker}(\mathrm{CH}^{(1)})\mathsf{p}}_{\boldsymbol{\mathcal{B}}^{(1)}, \varphi}$ is compatible with the $0$-source.}

Let $s$ be a sort in $S$ and let us consider the $0$-source operation symbol $\mathrm{sc}_{s}^{0}$ in $\Sigma^{\boldsymbol{\mathcal{A}}^{(1)}}_{s,s}$. Let $x$ be an element in $\T_{\boldsymbol{\mathcal{E}}^{\boldsymbol{\mathcal{B}}^{(1)}}}(\mathbf{Pth}_{\boldsymbol{\mathcal{B}}^{(1)}})_{\varphi(s)}$.

The following chain of equalities holds
\allowdisplaybreaks
\begin{align*}
\mathrm{pr}^{\mathrm{Ker}(\mathrm{CH}^{(1)})\mathsf{p}}_{\boldsymbol{\mathcal{B}}^{(1)}, \varphi(s)} \left(
\mathrm{sc}_{s}^{0\mathbf{T}_{\boldsymbol{\mathcal{E}}^{\boldsymbol{\mathcal{B}}^{(1)}}}^{\mathbf{f}^{(1)}}(\mathbf{Pth}_{\boldsymbol{\mathcal{B}}^{(1)}})}\left(
x
\right)
\right)
&=
\mathrm{pr}^{\mathrm{Ker}(\mathrm{CH}^{(1)})\mathsf{p}}_{\boldsymbol{\mathcal{B}}^{(1)}, \varphi(s)} \left(
\mathrm{sc}_{\varphi(s)}^{0\mathbf{T}_{\boldsymbol{\mathcal{E}}^{\boldsymbol{\mathcal{B}}^{(1)}}}(\mathbf{Pth}_{\boldsymbol{\mathcal{B}}^{(1)}})}\left(
x
\right)
\right)
\tag{1}
\\
&=
\mathrm{sc}_{\varphi(s)}^{0[\mathbf{Pth}_{\boldsymbol{\mathcal{B}}^{(1)}}]}\left(
\mathrm{pr}^{\mathrm{Ker}(\mathrm{CH}^{(1)})\mathsf{p}}_{\boldsymbol{\mathcal{B}}^{(1)}, \varphi(s)} \left(
x
\right)
\right)
\tag{2}
\\
&=
\mathrm{sc}_{s}^{0[\mathbf{Pth}_{\boldsymbol{\mathcal{B}}^{(1)}}^{\mathbf{f}^{(1)}}]}\left(
\mathrm{pr}^{\mathrm{Ker}(\mathrm{CH}^{(1)})\mathsf{p}}_{\boldsymbol{\mathcal{B}}^{(1)}, \varphi(s)} \left(
x
\right)
\right).
\tag{3}
\end{align*}
The first equality unravels the interpretation of the operation symbol $\mathrm{sc}_{s}^{0}$ in the partial $\Sigma^{\boldsymbol{\mathcal{A}}^{(1)}}$-algebra $\mathbf{T}_{\boldsymbol{\mathcal{E}}^{\boldsymbol{\mathcal{B}}^{(1)}}}^{\mathbf{f}^{(1)}}(\mathbf{Pth}_{\boldsymbol{\mathcal{B}}^{(1)}})$, introduced in Proposition~\ref{PFreeBCatAlg};
the second equality follows from the fact that, according to Proposition~\ref{CVarPr}, $\mathrm{pr}^{\mathrm{Ker}(\mathrm{CH}^{(1)}_{\boldsymbol{\mathcal{A}}^{(1)}})\mathsf{p}}_{\boldsymbol{\mathcal{B}}^{(1)}}$ is a $\Lambda^{\boldsymbol{\mathcal{B}}^{(1)}}$-homomorphism from $\mathbf{T}_{\boldsymbol{\mathcal{E}}^{\boldsymbol{\mathcal{B}}^{(1)}}}(\mathbf{Pth}_{\boldsymbol{\mathcal{B}}^{(1)}})$ to $[\mathbf{Pth}_{\boldsymbol{\mathcal{B}}^{(1)}}]$;
finally, the last equality recovers the interpretation of the operation symbol $\mathrm{sc}_{s}^{0}$ in the partial $\Sigma^{\boldsymbol{\mathcal{A}}^{(1)}}$-algebra $[\mathbf{Pth}_{\boldsymbol{\mathcal{B}}^{(1)}}^{\mathbf{f}^{(1)}}]$, introduced in Proposition~\ref{PQPthBCatAlg}.

Hence, $\mathrm{pr}^{\mathrm{Ker}(\mathrm{CH}^{(1)})\mathsf{p}}_{\boldsymbol{\mathcal{B}}^{(1)}, \varphi}$ is compatible with the $0$-source operation.

{\sffamily The mapping $\mathrm{pr}^{\mathrm{Ker}(\mathrm{CH}^{(1)})\mathsf{p}}_{\boldsymbol{\mathcal{B}}^{(1)}, \varphi}$ is compatible with the $0$-target.}

Let $s$ be a sort in $S$ and let us consider the $0$-target operation symbol $\mathrm{tg}_{s}^{0}$ in $\Sigma^{\boldsymbol{\mathcal{A}}^{(1)}}_{s,s}$. Let $x$ be an element in $\T_{\boldsymbol{\mathcal{E}}^{\boldsymbol{\mathcal{B}}^{(1)}}}(\mathbf{Pth}_{\boldsymbol{\mathcal{B}}^{(1)}})_{\varphi(s)}$, then the following equality holds
$$
\mathrm{pr}^{\mathrm{Ker}(\mathrm{CH}^{(1)})\mathsf{p}}_{\boldsymbol{\mathcal{B}}^{(1)}, \varphi(s)} \left(
\mathrm{tg}_{s}^{0\mathbf{T}_{\boldsymbol{\mathcal{E}}^{\boldsymbol{\mathcal{B}}^{(1)}}}^{\mathbf{f}^{(1)}}(\mathbf{Pth}_{\boldsymbol{\mathcal{B}}^{(1)}})}\left(
x
\right)
\right)
=
\mathrm{tg}_{s}^{0[\mathbf{Pth}_{\boldsymbol{\mathcal{B}}^{(1)}}^{\mathbf{f}^{(1)}}]}\left(
\mathrm{pr}^{\mathrm{Ker}(\mathrm{CH}^{(1)})\mathsf{p}}_{\boldsymbol{\mathcal{B}}^{(1)}, \varphi(s)} \left(
x
\right)
\right).
$$

The proof of this case is identical to that of the $0$-source.

Hence, $\mathrm{pr}^{\mathrm{Ker}(\mathrm{CH}^{(1)})\mathsf{p}}_{\boldsymbol{\mathcal{B}}^{(1)}, \varphi}$ is compatible with the $0$-target operation.

{\sffamily The mapping $\mathrm{pr}^{\mathrm{Ker}(\mathrm{CH}^{(1)})\mathsf{p}}_{\boldsymbol{\mathcal{B}}^{(1)}, \varphi}$ is compatible with the $0$-composition.}

Let $s$ be a sort in $S$ and let us consider the $0$-composition operation symbol $\circ_{s}^{0}$ in $\Sigma_{ss,s}^{\boldsymbol{\mathcal{A}}^{(1)}}$. Let $x$ and $y$ be two elements in $\T_{\boldsymbol{\mathcal{E}}^{\boldsymbol{\mathcal{B}}^{(1)}}}(\mathbf{Pth}_{\boldsymbol{\mathcal{B}}^{(1)}})_{\varphi(s)}$ such that
$$
\mathrm{sc}_{\varphi(s)}^{(0,1)\mathbf{T}_{\boldsymbol{\mathcal{E}}^{\boldsymbol{\mathcal{B}}^{(1)}}}^{\mathbf{f}^{(1)}}(\mathbf{Pth}_{\boldsymbol{\mathcal{B}}^{(1)}})}\left(
y
\right)
=
\mathrm{tg}_{\varphi(s)}^{(0,1)\mathbf{T}_{\boldsymbol{\mathcal{E}}^{\boldsymbol{\mathcal{B}}^{(1)}}}^{\mathbf{f}^{(1)}}(\mathbf{Pth}_{\boldsymbol{\mathcal{B}}^{(1)}})}\left(
x
\right).
$$
Then the following equality holds
$$
\mathrm{pr}^{\mathrm{Ker}(\mathrm{CH}^{(1)})\mathsf{p}}_{\boldsymbol{\mathcal{B}}^{(1)}, \varphi(s)} \left(
y
\circ_{s}^{0\mathbf{T}_{\boldsymbol{\mathcal{E}}^{\boldsymbol{\mathcal{B}}^{(1)}}}^{\mathbf{f}^{(1)}}(\mathbf{Pth}_{\boldsymbol{\mathcal{B}}^{(1)}})}\left(
x
\right)
\right)
=
\mathrm{pr}^{\mathrm{Ker}(\mathrm{CH}^{(1)})\mathsf{p}}_{\boldsymbol{\mathcal{B}}^{(1)}, \varphi(s)} \left(
y
\right)
\circ_{s}^{0[\mathbf{Pth}_{\boldsymbol{\mathcal{B}}^{(1)}}^{\mathbf{f}^{(1)}}]}
\mathrm{pr}^{\mathrm{Ker}(\mathrm{CH}^{(1)})\mathsf{p}}_{\boldsymbol{\mathcal{B}}^{(1)}, \varphi(s)} \left(
x
\right).
$$

The proof of this case is identical to that of the $0$-source.

Hence, $\mathrm{pr}^{\mathrm{Ker}(\mathrm{CH}^{(1)})\mathsf{p}}_{\boldsymbol{\mathcal{B}}^{(1)}, \varphi}$ is compatible with the $0$-composition operation.

This completes the proof.
\end{proof}

\begin{theorem}\label{TPthBFreeB}
The partial $\Sigma^{\boldsymbol{\mathcal{A}}^{(1)}}$-algebras $\mathbf{T}_{\boldsymbol{\mathcal{E}}^{\boldsymbol{\mathcal{B}}^{(1)}}}^{\mathbf{f}^{(1)}}(\mathbf{Pth}_{\boldsymbol{\mathcal{B}}^{(1)}})$ and $[\mathbf{Pth}_{\boldsymbol{\mathcal{B}}^{(1)}}^{\mathbf{f}^{(1)}}]$ are isomorphic.
\end{theorem}

\begin{proof}
Note that it follows from Theorem~\ref{TPthFree} that $(\mathrm{pr}^{\equiv^{[1]}}_{\boldsymbol{\mathcal{B}}^{(1)}} \circ \mathrm{ip}^{(1,Y)@}_{\boldsymbol{\mathcal{B}}^{(1)}} \circ \mathrm{CH}^{(1)}_{\boldsymbol{\mathcal{B}}^{(1)}})^{\natural}$ and $\mathrm{pr}^{\mathrm{Ker}(\mathrm{CH}^{(1)})\mathsf{p}}_{\boldsymbol{\mathcal{B}}^{(1)}}$ are a pair of inverse $\Lambda^{\boldsymbol{\mathcal{B}}^{(1)}}$-isomorphisms.

Thus, for every sort $s$ in $S$, the following chain of equalities holds
\begin{flushleft}
$
(\mathrm{pr}^{\equiv^{[1]}}_{\boldsymbol{\mathcal{B}}^{(1)}} \circ \mathrm{ip}^{(1,Y)@}_{\boldsymbol{\mathcal{B}}^{(1)}} \circ \mathrm{CH}^{(1)}_{\boldsymbol{\mathcal{B}}^{(1)}})^{\natural}_{\varphi(s)}
\circ
\mathrm{pr}^{\mathrm{Ker}(\mathrm{CH}^{(1)})\mathsf{p}}_{\boldsymbol{\mathcal{B}}^{(1)}, \varphi(s)}
$
\allowdisplaybreaks
\begin{align*}
&=
\left(
(\mathrm{pr}^{\equiv^{[1]}}_{\boldsymbol{\mathcal{B}}^{(1)}} \circ \mathrm{ip}^{(1,Y)@}_{\boldsymbol{\mathcal{B}}^{(1)}} \circ \mathrm{CH}^{(1)}_{\boldsymbol{\mathcal{B}}^{(1)}})^{\natural}
\circ
\mathrm{pr}^{\mathrm{Ker}(\mathrm{CH}^{(1)}_{\boldsymbol{\mathcal{A}}^{(1)}})\mathsf{p}}_{\boldsymbol{\mathcal{B}}^{(1)}}
\right)_{\varphi(s)}
\tag{1}
\\
&=
\left(
\T_{\boldsymbol{\mathcal{E}}^{\boldsymbol{\mathcal{B}}^{(1)}}}(\mathbf{Pth}_{\boldsymbol{\mathcal{B}}^{(1)}})
\right)_{\varphi(s)}
\tag{2}
\end{align*}
\end{flushleft}

Likewise, for every sort $s$ in $S$, the following chain of equalities holds
\begin{flushleft}
$
\mathrm{pr}^{\mathrm{Ker}(\mathrm{CH}^{(1)})\mathsf{p}}_{\boldsymbol{\mathcal{B}}^{(1)}, \varphi(s)}
\circ
(\mathrm{pr}^{\equiv^{[1]}}_{\boldsymbol{\mathcal{B}}^{(1)}} \circ \mathrm{ip}^{(1,Y)@}_{\boldsymbol{\mathcal{B}}^{(1)}} \circ \mathrm{CH}^{(1)}_{\boldsymbol{\mathcal{B}}^{(1)}})^{\natural}_{\varphi(s)}
$
\begin{align*}
&=
\left(
\mathrm{pr}^{\mathrm{Ker}(\mathrm{CH}^{(1)}_{\boldsymbol{\mathcal{A}}^{(1)}})\mathsf{p}}_{\boldsymbol{\mathcal{B}}^{(1)}}
\circ
(\mathrm{pr}^{\equiv^{[1]}}_{\boldsymbol{\mathcal{B}}^{(1)}} \circ \mathrm{ip}^{(1,Y)@}_{\boldsymbol{\mathcal{B}}^{(1)}} \circ \mathrm{CH}^{(1)}_{\boldsymbol{\mathcal{B}}^{(1)}})^{\natural}
\right)_{\varphi(s)}
\tag{1}
\\
&=
\left(
[\mathrm{Pth}_{\boldsymbol{\mathcal{B}}^{(1)}}]
\right)_{\varphi(s)}
\tag{2}
\end{align*}
\end{flushleft}
\end{proof}

\begin{corollary}\label{TPTBFreeB}
The partial $\Sigma^{\boldsymbol{\mathcal{A}}^{(1)}}$-algebras $\mathbf{T}_{\boldsymbol{\mathcal{E}}^{\boldsymbol{\mathcal{B}}^{(1)}}}^{\mathbf{f}^{(1)}}(\mathbf{Pth}_{\boldsymbol{\mathcal{B}}^{(1)}})$ and $[\mathbf{PT}_{\boldsymbol{\mathcal{B}}^{(1)}}^{\mathbf{f}^{(1)}}]$ are isomorphic.
\end{corollary}			
\chapter{First-order quotient path extension mapping}\label{S3C}

The aim of this chapter is to define, given a first-order morphism $\mathbf{f}^{(1)}$ and its path extension $f^{(1)\flat}$, a mapping from the $S$-sorted set $[ \mathrm{Pth}_{\boldsymbol{\mathcal{A}}^{(1)}} ]$ to the $S$-sorted set $[ \mathrm{Pth}_{\boldsymbol{\mathcal{B}}^{(1)}} ]_{\varphi}$ that we will call the first-order quotient path-extension mapping denoted by $f^{[ 1 ] @}$ which is a $\Sigma^{\boldsymbol{\mathcal{A}}^{(1)}}$-homomorphism from $[\mathbf{Pth}_{\boldsymbol{\mathcal{A}}^{(1)}}]$ to $[\mathbf{Pth}_{\boldsymbol{\mathcal{B}}^{(1)}}^{\mathbf{f}^{(1)}}]$.

In order to construct the desired mapping, we begin by proving that the partial $\Sigma^{\boldsymbol{\mathcal{A}}^{(1)}}$-algebra $[ \mathbf{Pth}_{\boldsymbol{\mathcal{B}}^{(1)}}^{\mathbf{f}^{(1)}} ]$ satisfies axiom A8 defining the QE-variety $\mathcal{V}(\boldsymbol{\mathcal{E}}^{\boldsymbol{\mathcal{A}}^{(1)}})$, introduced in Definition~\ref{DVar}.

\begin{proposition}
\label{PQPthBVarA8}
Let $(\mathbf{s}, s)$ an element of $S^{\star} \times S$, $\sigma$ an operation symbol in $\Sigma_{\mathbf{s}, s}$ and $( [\mathfrak{P}_{j}]_{\varphi(s_{j})} )_{j\in\bb{\mathbf{s}}}$ and $( [\mathfrak{Q}_{j}]_{\varphi(s_{j})} )_{j\in\bb{\mathbf{s}}}$ be two family of path classes in $[\mathrm{Pth}_{\boldsymbol{\mathcal{B}}^{(1)}}]_{\varphi^{\star}(\mathbf{s})}$ such that, for every $j \in \bb{\mathbf{s}}$, 
$$
\mathrm{sc}_{s_{j}}^{0[\mathbf{Pth}_{\boldsymbol{\mathcal{B}}^{(1)}}^{\mathbf{f}^{(1)}}]} \left(
\left[
\mathfrak{Q}_{j}
\right]_{s_{j}}
\right)
=
\mathrm{tg}_{s_{j}}^{0[\mathbf{Pth}_{\boldsymbol{\mathcal{B}}^{(1)}}^{\mathbf{f}^{(1)}}]} \left(
\left[
\mathfrak{P}_{j}
\right]_{s_{j}}
\right).
$$
Then the following equality holds
\begin{multline*}
\sigma^{[\mathbf{Pth}_{\boldsymbol{\mathcal{B}}^{(1)}}^{\mathbf{f}^{(1)}}]} \left(
\left(
\left[
\mathfrak{Q}_{j}
\right]_{s_{j}}
\circ_{s_{j}}^{0[\mathbf{Pth}_{\boldsymbol{\mathcal{B}}^{(1)}}^{\mathbf{f}^{(1)}}]}
\left[
\mathfrak{P}_{j}
\right]_{s_{j}}
\right)_{j \in \bb{\mathbf{s}}}
\right)
\\
=
\sigma^{[\mathbf{Pth}_{\boldsymbol{\mathcal{B}}^{(1)}}^{\mathbf{f}^{(1)}}]} \left(
\left(
\left[
\mathfrak{Q}_{j}
\right]_{s_{j}}
\right)_{j \in \bb{\mathbf{s}}}
\right)
\circ_{s}^{0[\mathbf{Pth}_{\boldsymbol{\mathcal{B}}^{(1)}}^{\mathbf{f}^{(1)}}]}
\sigma^{[\mathbf{Pth}_{\boldsymbol{\mathcal{B}}^{(1)}}^{\mathbf{f}^{(1)}}]} \left(
\left(
\left[
\mathfrak{P}_{j}
\right]_{s_{j}}
\right)_{j \in \bb{\mathbf{s}}}
\right).
\end{multline*}
\end{proposition}

\begin{proof}
Let $(\mathbf{s}, s)$ an element of $S^{\star} \times S$, $\sigma$ an operation symbol in $\Sigma_{\mathbf{s}, s}$ and $( [\mathfrak{P}_{j}]_{\varphi(s_{j})} )_{j\in\bb{\mathbf{s}}}$ and $( [\mathfrak{Q}_{j}]_{\varphi(s_{j})} )_{j\in\bb{\mathbf{s}}}$ be two family of path classes in $[\mathrm{Pth}_{\boldsymbol{\mathcal{B}}^{(1)}}]_{\varphi^{\star}(\mathbf{s})}$ such that, for every $j \in \bb{\mathbf{s}}$, 
$$
\mathrm{sc}_{s_{j}}^{0[\mathbf{Pth}_{\boldsymbol{\mathcal{B}}^{(1)}}^{\mathbf{f}^{(1)}}]} \left(
\left[
\mathfrak{Q}_{j}
\right]_{s_{j}}
\right)
=
\mathrm{tg}_{s_{j}}^{0[\mathbf{Pth}_{\boldsymbol{\mathcal{B}}^{(1)}}^{\mathbf{f}^{(1)}}]} \left(
\left[
\mathfrak{P}_{j}
\right]_{s_{j}}
\right).
$$
Unraveling the definition of the interpretation of the operation symbols $\sigma$ and $\circ_{s}^{0}$ in the $\Sigma^{\boldsymbol{\mathcal{A}}^{(1)}}$-algebras $[\mathbf{Pth}_{\boldsymbol{\mathcal{B}}^{(1)}}^{\mathbf{f}^{(1)}}]$ and $\mathbf{Pth}_{\boldsymbol{\mathcal{B}}^{(1)}}^{\mathbf{f}^{(1)}}$ introduced in Propositions~\ref{PQPthBCatAlg} and \ref{PPthBCatAlg}, respectively, the desired equality is equivalent to
\begin{multline*}
\left[
\sigma^{\mathbf{c}_{\mathfrak{d}}^{\ast}(\mathbf{Pth}_{\boldsymbol{\mathcal{B}}^{(1)}}^{(0,1)})} \left(
\left(
\mathfrak{Q}_{j}
\circ_{\varphi(s_{j})}^{0\mathbf{Pth}_{\boldsymbol{\mathcal{B}}^{(1)}}}
\mathfrak{P}_{j}
\right)_{j\in\bb{\mathbf{s}}}
\right)
\right]_{\varphi(s)}
\\
=
\left[
\sigma^{\mathbf{c}_{\mathfrak{d}}^{\ast}(\mathbf{Pth}_{\boldsymbol{\mathcal{B}}^{(1)}}^{(0,1)})} \left(
\left(
\mathfrak{Q}_{j}
\right)_{j \in \bb{\mathbf{s}}}
\right)
\circ_{\varphi(s)}^{0\mathbf{Pth}_{\boldsymbol{\mathcal{B}}^{(1)}}}
\sigma^{\mathbf{c}_{\mathfrak{d}}^{\ast}(\mathbf{Pth}_{\boldsymbol{\mathcal{B}}^{(1)}}^{(0,1)})} \left(
\left(
\mathfrak{P}_{j}
\right)_{j \in \bb{\mathbf{s}}}
\right)
\right]_{\varphi(s)}.
\tag{\dag}
\end{multline*}

Let us recall that, the interpretation of the operation symbol $\sigma$ in the $\Sigma$-algebra $\mathbf{c}_{\mathfrak{d}}^{\ast}(\mathbf{Pth}_{\boldsymbol{\mathcal{B}}^{(1)}}^{(0,1)})$ is given by the derived operation on $\mathbf{Pth}_{\boldsymbol{\mathcal{B}}^{(1)}}^{(0,1)}$ determined by $c_{\mathbf{s}, s}(\sigma)$ which following Remark~\ref{Rdop}, we denote by $c(\sigma)^{\mathbf{Pth}_{\boldsymbol{\mathcal{B}}^{(1)}}^{(0,1)}}$. For the details about derived operations, see Definition~\ref{wdop}.

Thus, we prove that, for every $(\mathbf{t}, t)$ in $T^{\star} \times T$, every $P$ in $T_{\Lambda}(\vs\mathbf{t})_{t}$ and every pair of families $(\mathfrak{P}_{j})_{j \in \bb{\mathbf{t}}}$ and $(\mathfrak{Q}_{j})_{j \in \bb{\mathbf{t}}}$ in $\mathrm{Pth}_{\boldsymbol{\mathcal{B}}^{(1)}, \mathbf{t}}$,
\begin{multline*}
\left[
P^{\mathbf{Pth}_{\boldsymbol{\mathcal{B}}^{(1)}}^{(0,1)}} \left(
\left(
\mathfrak{Q}_{j}
\circ_{t_{j}}^{0\mathbf{Pth}_{\boldsymbol{\mathcal{B}}^{(1)}}}
\mathfrak{P}_{j}
\right)_{j\in\bb{\mathbf{t}}}
\right)
\right]_{t}
\\
=
\left[
P^{\mathbf{Pth}_{\boldsymbol{\mathcal{B}}^{(1)}}^{(0,1)}} \left(
\left(
\mathfrak{Q}_{j}
\right)_{j \in \bb{\mathbf{t}}}
\right)
\circ_{t}^{0\mathbf{Pth}_{\boldsymbol{\mathcal{B}}^{(1)}}}
P^{\mathbf{Pth}_{\boldsymbol{\mathcal{B}}^{(1)}}^{(0,1)}} \left(
\left(
\mathfrak{P}_{j}
\right)_{j \in \bb{\mathbf{t}}}
\right)
\right]_{t}.
\end{multline*}
Therefore, Equation~(\dag) follows as a consequence.

By algebraic induction on the complexity of $P$.

{\sffamily Base step of the induction}

If $P$ is a variable $v_{i}^{t}$, thus, $i \in \bb{t}$ and $t_{i} = t$, then following Remark~\ref{Rdop}
$$
(v_{i}^{t})^{\mathbf{Pth}_{\boldsymbol{\mathcal{B}}^{(1)}}^{(0,1)}}
=
\mathrm{d}_{\mathbf{t}, t}^{\mathbf{Pth}_{\boldsymbol{\mathcal{B}}^{(1)}}^{(0,1)}} (v_{i}^{t})
=
\mathrm{d}_{\mathbf{t}, t}^{\mathbf{Pth}_{\boldsymbol{\mathcal{B}}^{(1)}}^{(0,1)}} ( \eta^{\vs\mathbf{t}} (v_{i}^{t})) 
=
\mathrm{p}_{\mathbf{t}, t}^{\mathbf{Pth}_{\boldsymbol{\mathcal{B}}^{(1)}}^{(0,1)}} (v_{i}^{t})
=
\mathrm{pr}_{\mathbf{t}, i}^{\mathrm{Pth}_{\boldsymbol{\mathcal{B}}^{(1)}}}.
$$

The following chain of equalities holds
\allowdisplaybreaks
\begin{flushleft}
$
P^{\mathbf{Pth}_{\boldsymbol{\mathcal{B}}^{(1)}}^{(0,1)}} \left(
\left(
\mathfrak{Q}_{j}
\circ_{t_{j}}^{0\mathbf{Pth}_{\boldsymbol{\mathcal{B}}^{(1)}}}
\mathfrak{P}_{j}
\right)_{j\in\bb{\mathbf{t}}}
\right)
$
\allowdisplaybreaks
\begin{align*}
&=
(v_{i}^{t})^{\mathbf{Pth}_{\boldsymbol{\mathcal{B}}^{(1)}}^{(0,1)}} \left(
\left(
\mathfrak{Q}_{j}
\circ_{t_{j}}^{0\mathbf{Pth}_{\boldsymbol{\mathcal{B}}^{(1)}}}
\mathfrak{P}_{j}
\right)_{j\in\bb{\mathbf{t}}}
\right)
\tag{1}
\\
&=
\mathrm{pr}_{\mathbf{t}, i}^{\mathrm{Pth}_{\boldsymbol{\mathcal{B}}^{(1)}}} \left(
\left(
\mathfrak{Q}_{j}
\circ_{t_{j}}^{0\mathbf{Pth}_{\boldsymbol{\mathcal{B}}^{(1)}}}
\mathfrak{P}_{j}
\right)_{j\in\bb{\mathbf{t}}}
\right)
\tag{2}
\\
&=
\mathfrak{Q}_{i}
\circ_{t_{i}}^{0\mathbf{Pth}_{\boldsymbol{\mathcal{B}}^{(1)}}}
\mathfrak{P}_{i}
\tag{3}
\\
&=
\mathfrak{Q}_{i}
\circ_{t}^{0\mathbf{Pth}_{\boldsymbol{\mathcal{B}}^{(1)}}}
\mathfrak{P}_{i}
\tag{4}
\\
&=
\mathrm{pr}_{\mathbf{t}, i}^{\mathrm{Pth}_{\boldsymbol{\mathcal{B}}^{(1)}}} \left(
\left(
\mathfrak{Q}_{j}
\right)_{j \in \bb{\mathbf{t}}}
\right)
\circ_{t}^{0\mathbf{Pth}_{\boldsymbol{\mathcal{B}}^{(1)}}}
\mathrm{pr}_{\mathbf{t}, i}^{\mathrm{Pth}_{\boldsymbol{\mathcal{B}}^{(1)}}} \left(
\left(
\mathfrak{P}_{j}
\right)_{j\in\bb{\mathbf{t}}}
\right)
\tag{5}
\\
&=
(v_{i}^{t})^{\mathbf{Pth}_{\boldsymbol{\mathcal{B}}^{(1)}}^{(0,1)}} \left(
\left(
\mathfrak{Q}_{j}
\right)_{j \in \bb{\mathbf{t}}}
\right)
\circ_{t}^{0\mathbf{Pth}_{\boldsymbol{\mathcal{B}}^{(1)}}}
(v_{i}^{t})^{\mathbf{Pth}_{\boldsymbol{\mathcal{B}}^{(1)}}^{(0,1)}} \left(
\left(
\mathfrak{P}_{j}
\right)_{j\in\bb{\mathbf{t}}}
\right)
\tag{6}
\\
&=
P^{\mathbf{Pth}_{\boldsymbol{\mathcal{B}}^{(1)}}^{(0,1)}} \left(
\left(
\mathfrak{Q}_{j}
\right)_{j \in \bb{\mathbf{t}}}
\right)
\circ_{t}^{0\mathbf{Pth}_{\boldsymbol{\mathcal{B}}^{(1)}}}
P^{\mathbf{Pth}_{\boldsymbol{\mathcal{B}}^{(1)}}^{(0,1)}} \left(
\left(
\mathfrak{P}_{j}
\right)_{j\in\bb{\mathbf{t}}}
\right).
\tag{7}
\end{align*}
\end{flushleft}

The first equality unravels the definition of $P$;
the second equality follows from the fact that $ (v_{i}^{t})^{\mathbf{Pth}_{\boldsymbol{\mathcal{B}}^{(1)}}^{(0,1)}} = \mathrm{pr}_{i}^{\mathrm{Pth}_{\boldsymbol{\mathcal{B}}^{(1)}, \mathbf{t}}} $;
the third equality unravels the definition of the $\mathbf{t}$-ary, $i$-th canonical projection from $\mathrm{Pth}_{\boldsymbol{\mathcal{B}}^{(1)}, \mathbf{t}}$ to $\mathrm{Pth}_{\boldsymbol{\mathcal{B}}^{(1)}, t_{i}}$;
the fourth equality follows from the fact that $t_{i} = t$;
the fifth equality recovers the definition of the $\mathbf{t}$-ary, $i$-th canonical projection from $\mathrm{Pth}_{\boldsymbol{\mathcal{B}}^{(1)}, \mathbf{t}}$ to $\mathrm{Pth}_{\boldsymbol{\mathcal{B}}^{(1)}, t_{i}}$;
the sixth equality follows from the fact that $ (v_{i}^{t})^{\mathbf{Pth}_{\boldsymbol{\mathcal{B}}^{(1)}}^{(0,1)}} = \mathrm{pr}_{i}^{\mathrm{Pth}_{\boldsymbol{\mathcal{B}}^{(1)}, \mathbf{t}}} $;
finally, the last equality recovers the definition of $P$.

Therefore, its $\mathrm{Ker}(\mathrm{CH}^{(1)}_{\boldsymbol{\mathcal{B}}^{(1)}})$-classes coincide, i.e., equality (\dag) follows.

If $P$ is a constant symbol $\tau$ in $\Lambda_{\lambda, t}$, then the following chain of equalities holds
\allowdisplaybreaks
\begin{flushleft}
$
P^{\mathbf{Pth}_{\boldsymbol{\mathcal{B}}^{(1)}}^{(0,1)}} \left(
\left(
\mathfrak{Q}_{j}
\circ_{t_{j}}^{0\mathbf{Pth}_{\boldsymbol{\mathcal{B}}^{(1)}}}
\mathfrak{P}_{j}
\right)_{j\in\bb{\mathbf{t}}}
\right)
$
\begin{align*}
&=
\tau^{\mathbf{Pth}_{\boldsymbol{\mathcal{B}}^{(1)}}^{(0,1)}}
\tag{1}
\\
&=
\tau^{\mathbf{Pth}_{\boldsymbol{\mathcal{B}}^{(1)}}^{(0,1)}}
\circ_{t}^{0\mathbf{Pth}_{\boldsymbol{\mathcal{B}}^{(1)}}}
\tau^{\mathbf{Pth}_{\boldsymbol{\mathcal{B}}^{(1)}}^{(0,1)}}
\tag{2}
\\
&=
P^{\mathbf{Pth}_{\boldsymbol{\mathcal{B}}^{(1)}}^{(0,1)}} \left(
\left(
\mathfrak{Q}_{j}
\right)_{j \in \bb{\mathbf{t}}}
\right)
\circ_{t}^{0\mathbf{Pth}_{\boldsymbol{\mathcal{B}}^{(1)}}}
P^{\mathbf{Pth}_{\boldsymbol{\mathcal{B}}^{(1)}}^{(0,1)}} \left(
\left(
\mathfrak{P}_{j}
\right)_{j\in\bb{\mathbf{t}}}
\right).
\tag{3}
\end{align*}
\end{flushleft}

The first equality unravels the definition of $P$;
the second equality follows from the fact that, according to Proposition~\ref{PPthCatAlg}, the interpretation of the $\tau$ operation in the many-sorted partial $\Lambda^{\boldsymbol{\mathcal{B}}^{(1)}}$-algebra $\mathbf{Pth}_{\boldsymbol{\mathcal{B}}^{(1)}}$ is given by the $(1,0)$-identity path on $[\tau^{\mathbf{T}_{\Lambda}(Y)}]$, i.e., $\tau^{\mathbf{Pth}_{\boldsymbol{\mathcal{B}}^{(1)}}} = \mathrm{ip}_{t}^{(1,0)\sharp} \left([\tau^{\mathbf{T}_{\Lambda}(Y)}]\right)$ and from the fact that, according to Proposition~\ref{PPthComp}, the $(1,0)$-identity paths are idempotent for the $0$-composition;
finally, the last equality recovers the definition of $P$.

Therefore, its $\mathrm{Ker}(\mathrm{CH}^{(1)}_{\boldsymbol{\mathcal{B}}^{(1)}})$-classes coincide, i.e., equality (\dag) follows.

{\sffamily Inductive step of the induction}

If $P$ is a term $\tau \left(\left(Q_{i}\right)_{i \in \bb{\mathbf{r}}}\right)$ for $\tau$ an operation symbol in $\lambda_{\mathbf{r}, t}$ and $(Q_{i})_{i \in \bb{\mathbf{r}}}$ a family of terms in $\T_{\Lambda}(\vs\mathbf{t})_{\mathbf{r}}$, then following Remark~\ref{Rdop}
\begin{align*}
P^{\mathbf{Pth}_{\boldsymbol{\mathcal{B}}^{(1)}}^{(0,1)}}
&=
\mathrm{d}_{\mathbf{t}, t}^{\mathbf{Pth}_{\boldsymbol{\mathcal{B}}^{(1)}}^{(0,1)}} (P)
\\
&=
\mathrm{d}_{\mathbf{t}, t}^{\mathbf{Pth}_{\boldsymbol{\mathcal{B}}^{(1)}}^{(0,1)}} \left(\tau \left(\left(Q_{i}\right)_{i \in \bb{\mathbf{r}}}\right)\right) 
\\
&=
\tau^{\mathbf{O}_{\mathbf{t}}\left(\mathbf{Pth}_{\boldsymbol{\mathcal{B}}^{(1)}}^{(0,1)}\right)} \left(
\left(
\mathrm{d}_{\mathbf{t}, r_{i}}^{\mathbf{Pth}_{\boldsymbol{\mathcal{B}}^{(1)}}^{(0,1)}}
\left(
Q_{i}
\right)
\right)_{i \in \bb{\mathbf{r}}}
\right)
\\
&=
\tau^{\mathbf{O}_{\mathbf{t}}\left(\mathbf{Pth}_{\boldsymbol{\mathcal{B}}^{(1)}}^{(0,1)}\right)} \left(
\left(
Q_{i}^{\mathbf{Pth}_{\boldsymbol{\mathcal{B}}^{(1)}}^{(0,1)}}
\right)_{i \in \bb{\mathbf{r}}}
\right).
\end{align*}

The following chain of equalities holds
\allowdisplaybreaks
\begin{flushleft}
$
\left[
P^{\mathbf{Pth}_{\boldsymbol{\mathcal{B}}^{(1)}}^{(0,1)}} \left(
\left(
\mathfrak{Q}_{j}
\circ_{t_{j}}^{0\mathbf{Pth}_{\boldsymbol{\mathcal{B}}^{(1)}}}
\mathfrak{P}_{j}
\right)_{j\in\bb{\mathbf{t}}}
\right)
\right]_{t}
$
\begin{align*}
&=
\left[
\tau^{\mathbf{O}_{\mathbf{t}}\left(\mathbf{Pth}_{\boldsymbol{\mathcal{B}}^{(1)}}^{(0,1)}\right)} \left(
\left(
Q_{i}^{\mathbf{Pth}_{\boldsymbol{\mathcal{B}}^{(1)}}^{(0,1)}}
\right)_{i \in \bb{\mathbf{r}}}
\right)
\left(
\left(
\mathfrak{Q}_{j}
\circ_{t_{j}}^{0\mathbf{Pth}_{\boldsymbol{\mathcal{B}}^{(1)}}}
\mathfrak{P}_{j}
\right)_{j\in\bb{\mathbf{t}}}
\right)
\right]_{t}
\tag{1}
\\
&=
\left[
\tau^{\mathbf{Pth}_{\boldsymbol{\mathcal{B}}^{(1)}}^{(0,1)}} \left(
\left(
Q_{i}^{\mathbf{Pth}_{\boldsymbol{\mathcal{B}}^{(1)}}^{(0,1)}}
\left(
\left(
\mathfrak{Q}_{j}
\circ_{t_{j}}^{0\mathbf{Pth}_{\boldsymbol{\mathcal{B}}^{(1)}}}
\mathfrak{P}_{j}
\right)_{j\in\bb{\mathbf{t}}}
\right)
\right)_{i \in \bb{\mathbf{r}}}
\right)
\right]_{t}
\tag{2}
\\
&=
\tau^{\left[\mathbf{Pth}_{\boldsymbol{\mathcal{B}}^{(1)}}^{(0,1)}\right]} \left(
\left(
\left[
Q_{i}^{\mathbf{Pth}_{\boldsymbol{\mathcal{B}}^{(1)}}^{(0,1)}}
\left(
\left(
\mathfrak{Q}_{j}
\circ_{t_{j}}^{0\mathbf{Pth}_{\boldsymbol{\mathcal{B}}^{(1)}}}
\mathfrak{P}_{j}
\right)_{j\in\bb{\mathbf{t}}}
\right)
\right]_{r_{i}}
\right)_{i\in\bb{\mathbf{r}}}
\right)
\tag{3}
\\
&=
\resizebox{0.86\textwidth}{!}{%
$
\tau^{\left[\mathbf{Pth}_{\boldsymbol{\mathcal{B}}^{(1)}}^{(0,1)}\right]} \left(
\left(
\left[
Q_{i}^{\mathbf{Pth}_{\boldsymbol{\mathcal{B}}^{(1)}}^{(0,1)}} \left(
\left(
\mathfrak{Q}_{j}
\right)_{j \in \bb{\mathbf{t}}}
\right)
\circ_{r_{i}}^{0\mathbf{Pth}_{\boldsymbol{\mathcal{B}}^{(1)}}}
Q_{i}^{\mathbf{Pth}_{\boldsymbol{\mathcal{B}}^{(1)}}^{(0,1)}} \left(
\left(
\mathfrak{P}_{j}
\right)_{j \in \bb{\mathbf{t}}}
\right)
\right]_{r_{i}}
\right)_{i\in\bb{\mathbf{r}}}
\right)
$}
\tag{4}
\\
&=
\resizebox{0.86\textwidth}{!}{%
$
\tau^{\left[\mathbf{Pth}_{\boldsymbol{\mathcal{B}}^{(1)}}^{(0,1)}\right]} \left(
\left(
\left[
Q_{i}^{\mathbf{Pth}_{\boldsymbol{\mathcal{B}}^{(1)}}^{(0,1)}} \left(
\left(
\mathfrak{Q}_{j}
\right)_{j \in \bb{\mathbf{t}}}
\right)
\right]_{r_{i}}
\circ_{r_{i}}^{0\left[\mathbf{Pth}_{\boldsymbol{\mathcal{B}}^{(1)}}\right]}
\left[
Q_{i}^{\mathbf{Pth}_{\boldsymbol{\mathcal{B}}^{(1)}}^{(0,1)}} \left(
\left(
\mathfrak{P}_{j}
\right)_{j \in \bb{\mathbf{t}}}
\right)
\right]_{r_{i}}
\right)_{i\in\bb{\mathbf{r}}}
\right)
$}
\tag{5}
\\
&=
\tau^{\left[\mathbf{Pth}_{\boldsymbol{\mathcal{B}}^{(1)}}^{(0,1)}\right]} \left(
\left(
\left[
Q_{i}^{\mathbf{Pth}_{\boldsymbol{\mathcal{B}}^{(1)}}^{(0,1)}} \left(
\left(
\mathfrak{Q}_{j}
\right)_{j \in \bb{\mathbf{t}}}
\right)
\right]_{r_{i}}
\right)_{i\in\bb{\mathbf{r}}}
\right)
\\
&\hspace{2cm}
\circ_{t}^{0\left[\mathbf{Pth}_{\boldsymbol{\mathcal{B}}^{(1)}}\right]}
\tau^{\left[\mathbf{Pth}_{\boldsymbol{\mathcal{B}}^{(1)}}^{(0,1)}\right]} \left(
\left(
\left[
Q_{i}^{\mathbf{Pth}_{\boldsymbol{\mathcal{B}}^{(1)}}^{(0,1)}} \left(
\left(
\mathfrak{P}_{j}
\right)_{j \in \bb{\mathbf{t}}}
\right)
\right]_{r_{i}}
\right)_{i\in\bb{\mathbf{r}}}
\right)
\tag{6}
\\
&=
\left[
\tau^{\mathbf{Pth}_{\boldsymbol{\mathcal{B}}^{(1)}}^{(0,1)}} \left(
\left(
Q_{i}^{\mathbf{Pth}_{\boldsymbol{\mathcal{B}}^{(1)}}^{(0,1)}} \left(
\left(
\mathfrak{Q}_{j}
\right)_{j \in \bb{\mathbf{t}}}
\right)
\right)_{i\in\bb{\mathbf{r}}}
\right)
\right]_{t}
\\
&\hspace{2.5cm}
\circ_{t}^{0\left[\mathbf{Pth}_{\boldsymbol{\mathcal{B}}^{(1)}}\right]}
\left[
\tau^{\mathbf{Pth}_{\boldsymbol{\mathcal{B}}^{(1)}}^{(0,1)}} \left(
\left(
Q_{i}^{\mathbf{Pth}_{\boldsymbol{\mathcal{B}}^{(1)}}^{(0,1)}} \left(
\left(
\mathfrak{P}_{j}
\right)_{j \in \bb{\mathbf{t}}}
\right)
\right)_{i\in\bb{\mathbf{r}}}
\right)
\right]_{t}
\tag{7}
\\
&=
\left[
\tau^{\mathbf{Pth}_{\boldsymbol{\mathcal{B}}^{(1)}}^{(0,1)}} \left(
\left(
Q_{i}^{\mathbf{Pth}_{\boldsymbol{\mathcal{B}}^{(1)}}^{(0,1)}} \left(
\left(
\mathfrak{Q}_{j}
\right)_{j \in \bb{\mathbf{t}}}
\right)
\right)_{i\in\bb{\mathbf{r}}}
\right)
\right.
\\
&\hspace{3.1cm}
\left.
\circ_{t}^{0\mathbf{Pth}_{\boldsymbol{\mathcal{B}}^{(1)}}}
\tau^{\mathbf{Pth}_{\boldsymbol{\mathcal{B}}^{(1)}}^{(0,1)}} \left(
\left(
Q_{i}^{\mathbf{Pth}_{\boldsymbol{\mathcal{B}}^{(1)}}^{(0,1)}} \left(
\left(
\mathfrak{P}_{j}
\right)_{j \in \bb{\mathbf{t}}}
\right)
\right)_{i\in\bb{\mathbf{r}}}
\right)
\right]_{t}
\tag{8}
\\
&=
\left[
\tau^{\mathbf{O}_{\mathbf{t}}\left(\mathbf{Pth}_{\boldsymbol{\mathcal{B}}^{(1)}}^{(0,1)}\right)} \left(
\left(
Q_{i}^{\mathbf{Pth}_{\boldsymbol{\mathcal{B}}^{(1)}}^{(0,1)}}
\right)_{i \in \bb{\mathbf{r}}}
\right)
\left(
\left(
\mathfrak{Q}_{j}
\right)_{j\in\bb{\mathbf{t}}}
\right)
\right.
\\
&\hspace{2.5cm}
\left.
\circ_{t}^{0\mathbf{Pth}_{\boldsymbol{\mathcal{B}}^{(1)}}}
\tau^{\mathbf{O}_{\mathbf{t}}\left(\mathbf{Pth}_{\boldsymbol{\mathcal{B}}^{(1)}}^{(0,1)}\right)} \left(
\left(
Q_{i}^{\mathbf{Pth}_{\boldsymbol{\mathcal{B}}^{(1)}}^{(0,1)}}
\right)_{i \in \bb{\mathbf{r}}}
\right)
\left(
\left(
\mathfrak{P}_{j}
\right)_{j\in\bb{\mathbf{t}}}
\right)
\right]_{t}
\tag{9}
\\
&=
\left[
P^{\mathbf{Pth}_{\boldsymbol{\mathcal{B}}^{(1)}}^{(0,1)}}
\left(
\left(
\mathfrak{Q}_{j}
\right)_{j\in\bb{\mathbf{t}}}
\right)
\circ_{t}^{0\mathbf{Pth}_{\boldsymbol{\mathcal{B}}^{(1)}}}
P^{\mathbf{Pth}_{\boldsymbol{\mathcal{B}}^{(1)}}^{(0,1)}}
\left(
\left(
\mathfrak{P}_{j}
\right)_{j\in\bb{\mathbf{t}}}
\right)
\right]_{t}.
\tag{10}
\end{align*}
\end{flushleft}

The first equality unravels the definition of $P$;
the second equality unravels the definition of the interpretation of the operation symbol $\tau$ in the product $\Lambda$-algebra $\mathbf{O}_{\mathbf{t}}(\mathbf{Pth}_{\boldsymbol{\mathcal{B}}^{(1)}}^{(0,1)})$ introduced in Definition~\ref{wdop};
the third equality follows from the definition of the interpretation of the operation symbol $\tau$ in the quotient $\Lambda$-algebra $[\mathbf{Pth}_{\boldsymbol{\mathcal{B}}^{(1)}}^{(0,1)}]$ introduced in Proposition~\ref{PCHCatAlg};
the fourth equality follows by the inductive hypothesis;
the fifth equality follows from the definition of the interpretation of the operation symbol $\circ_{r_{i}}^{0}$ in the quotient $\Lambda^{\boldsymbol{\mathcal{B}}^{(1)}}$-algebra $[\mathbf{Pth}_{\boldsymbol{\mathcal{B}}^{(1)}}]$ introduced in Proposition~\ref{PCHCatAlg};
the sixth equality follows from Proposition~\ref{PCHVarA8};
the seventh equality recovers the definition of the interpretation of the operation symbol $\tau$ in the quotient $\Lambda$-algebra $[\mathbf{Pth}_{\boldsymbol{\mathcal{B}}^{(1)}}^{(0,1)}]$ introduced in Proposition~\ref{PCHCatAlg};
the eighth equality recovers the definition of the interpretation of the operation symbol $\circ_{t}^{0}$ in the quotient $\Lambda^{\boldsymbol{\mathcal{B}}^{(1)}}$-algebra $[\mathbf{Pth}_{\boldsymbol{\mathcal{B}}^{(1)}}]$ introduced in Proposition~\ref{PCHCatAlg};
the ninth equality recovers the definition of the interpretation of the operation symbol $\tau$ in the product $\Lambda$-algebra $\mathbf{O}_{\mathbf{t}}(\mathbf{Pth}_{\boldsymbol{\mathcal{B}}^{(1)}}^{(0,1)})$ introduced in Definition~\ref{wdop};
finally, the last equality recovers the definition of $P$.
\end{proof}

We now work towards the existence of the first-order quotient path-extension mapping.

\begin{proposition}\label{PHomPthExtKer}
Let $\mathbf{f}^{(1)} = (\varphi,c,(f^{(i)})_{i\in 2})$ be a first-order morphism from $\boldsymbol{\mathcal{A}}^{(1)}$ to $\boldsymbol{\mathcal{B}}^{(1)}$. Then
$$
\mathrm{pr}_{\boldsymbol{\mathcal{B}}^{(1)}, \varphi}^{\mathrm{Ker}(\mathrm{CH}^{(1)})} \circ f^{(1)\flat}
\colon
\mathrm{Pth}_{\boldsymbol{\mathcal{A}}^{(1)}}
\mor
\left[\mathrm{Pth}_{\boldsymbol{\mathcal{B}}^{(1)}}\right]_{\varphi}
$$
is a $\Sigma^{\boldsymbol{\mathcal{A}}^{(1)}}$-homomorphism from $\mathbf{Pth}_{\boldsymbol{\mathcal{A}}^{(1)}}$ to $[\mathbf{Pth}_{\boldsymbol{\mathcal{B}}^{(1)}}^{\mathbf{f}^{(1)}}]_{\varphi}$ satisfying that
$$
\mathrm{Ker}(\mathrm{CH}^{(1)}_{\boldsymbol{\mathcal{A}}^{(1)}})
\subseteq
\mathrm{Ker}\left(\mathrm{pr}_{\boldsymbol{\mathcal{B}}^{(1)}, \varphi}^{\mathrm{Ker}(\mathrm{CH}^{(1)})} \circ f^{(1)\flat}\right).
$$
For further reference, we will call this composition mapping, which will become useful in further constructions, the \emph{first-order homomorphic path extension mapping} of $\mathbf{f}^{(1)}$.
\end{proposition}

\begin{proof}
We prove that $\mathrm{pr}_{\boldsymbol{\mathcal{B}}^{(1)}, \varphi}^{\mathrm{Ker}(\mathrm{CH}^{(1)})} \circ f^{(1)\flat}$ is compatible with every operation symbol in $\Sigma^{\boldsymbol{\mathcal{A}}^{(1)}}$.

{\sffamily The mapping $\mathrm{pr}_{\boldsymbol{\mathcal{B}}^{(1)}, \varphi}^{\mathrm{Ker}(\mathrm{CH}^{(1)})} \circ f^{(1)\flat}$ is a $\Sigma$-homomorphism.}


By Proposition~\ref{PPthExtHom}, $f^{(1)\flat}$ is a $\Sigma$-homomorphism from $\mathbf{Pth}_{\boldsymbol{\mathcal{A}}^{2}}^{(0,1)}$ to $\mathbf{Pth}_{\boldsymbol{\mathcal{B}}^{(1)}}^{\mathbf{f}^{(1)}(0,1)}$. Also, by Proposition~\ref{PKerCHBCatHom}, $\mathrm{pr}_{\boldsymbol{\mathcal{B}}^{(1)}, \varphi}^{\mathrm{Ker}(\mathrm{CH}^{(1)})}$ is a $\Sigma$-homomorphism from $\mathbf{Pth}_{\boldsymbol{\mathcal{B}}^{(1)}}^{\mathbf{f}^{(1)}(0,1)}$ to $[\mathbf{Pth}_{\boldsymbol{\mathcal{B}}^{(1)}}^{\mathbf{f}^{(1)}(0,1)}]$. Thus, the composition $\mathrm{pr}_{\boldsymbol{\mathcal{B}}^{(1)}, \varphi}^{\mathrm{Ker}(\mathrm{CH}^{(1)})} \circ f^{(1)\flat}$ is a $\Sigma$-homomorphism from $\mathbf{Pth}_{\boldsymbol{\mathcal{A}}^{(1)}}^{(0,1)}$ to $[\mathbf{Pth}_{\boldsymbol{\mathcal{B}}^{(1)}}^{\mathbf{f}^{(1)}(0,1)}]$.

Therefore, for every word $\mathbf{s}$ in $S^{\star}$, every sort $s$ in $S$, every operation symbol $\sigma \in \Sigma_{\mathbf{s}, s}$ and every family $(\mathfrak{P}_{j})_{j \in \bb{\mathbf{s}}}$ in $\mathrm{Pth}_{\boldsymbol{\mathcal{A}}^{(1)}, \mathbf{s}}$,
$$
\left[
f_{s}^{(1)\flat}\left(
\sigma^{\mathbf{Pth}_{\boldsymbol{\mathcal{A}}^{(1)}}}\left(\left(
\mathfrak{P}_{j}
\right)_{j\in\bb{\mathbf{s}}}\right)
\right)
\right]_{\varphi(s)}
=
\sigma^{[\mathbf{Pth}_{\boldsymbol{\mathcal{B}}^{(1)}}^{\mathbf{f}^{(1)}}]} \left(\left(
\left[
f_{s_{j}}^{(1)\flat}\left(
\mathfrak{P}_{j}
\right)
\right]_{\varphi(s_{j})}
\right)_{j\in\bb{\mathbf{s}}}\right).
$$

{\sffamily The mapping $\mathrm{pr}_{\boldsymbol{\mathcal{B}}^{(1)}, \varphi}^{\mathrm{Ker}(\mathrm{CH}^{(1)})} \circ f^{(1)\flat}$ is compatible with the rewrite rules.}

Let $s$ be a sort in $S$ and let $\mathfrak{p}$ be a rewrite rule in $\mathcal{A}_{s}^{(1)}$. We have to show that
$$
\left[
f_{s}^{(1)\flat}\left(
\mathfrak{p}^{\mathbf{Pth}_{\boldsymbol{\mathcal{A}}^{(1)}}}
\right)
\right]_{\varphi(s)}
=
\mathfrak{p}^{[\mathbf{Pth}_{\boldsymbol{\mathcal{B}}^{(1)}}^{\mathbf{f}^{(1)}}]}
$$

The following chain of equalities holds
\allowdisplaybreaks
\begin{align*}
\left[
f_{s}^{(1)\flat}\left(
\mathfrak{p}^{\mathbf{Pth}_{\boldsymbol{\mathcal{A}}^{(1)}}}
\right)
\right]_{\varphi(s)}
&=
\left[
f_{s}^{(1)\flat}\left(
\mathrm{ech}_{\boldsymbol{\mathcal{A}}^{(1)}, s}^{(1,\mathcal{A}^{(1)})}\left(
\mathfrak{p}
\right)
\right)
\right]_{\varphi(s)}
\tag{1}
\\
&=
\mathfrak{p}^{[\mathbf{Pth}_{\boldsymbol{\mathcal{B}}^{(1)}}^{\mathbf{f}^{(1)}}]}.
\tag{2}
\end{align*}

The first equality unravels the interpretation of the constant symbol $\mathfrak{p}$ in the partial $\Sigma^{\boldsymbol{\mathcal{A}}^{(1)}}$-algebra $\mathbf{Pth}_{\boldsymbol{\mathcal{A}}^{(1)}}$, introduced in Proposition~\ref{PDPthCatAlg};
the second equality recovers the interpretation of the constant symbol $\mathfrak{p}$ in the partial $\Sigma^{\boldsymbol{\mathcal{A}}^{(1)}}$-algebra $[\mathbf{Pth}_{\boldsymbol{\mathcal{B}}^{(1)}}^{\mathbf{f}^{(1)}}]$, introduced in Proposition~\ref{PQPthBCatAlg}.

Hence, $\mathrm{pr}_{\boldsymbol{\mathcal{B}}^{(1)}, \varphi}^{\mathrm{Ker}(\mathrm{CH}^{(1)})} \circ f^{(1)\flat}$ is compatible with the rewrite rules.

{\sffamily The mapping $\mathrm{pr}_{\boldsymbol{\mathcal{B}}^{(1)}, \varphi}^{\mathrm{Ker}(\mathrm{CH}^{(1)})} \circ f^{(1)\flat}$ is compatible with the $0$-source.}

Let $s$ be a sort in $S$ and let us consider the $0$-source operation symbol $\mathrm{sc}_{s}^{0}$ in $\Sigma^{\boldsymbol{\mathcal{A}}^{(1)}}_{s,s}$. Let $\mathfrak{P}$ be a path in $\mathrm{Pth}_{\boldsymbol{\mathcal{A}}^{(1)}, s}$.

The following chain of equalities holds
\allowdisplaybreaks
\begin{align*}
\left[
f_{s}^{(1)\flat} \left(
\mathrm{sc}_{s}^{0\mathbf{Pth}_{\boldsymbol{\mathcal{A}}^{(1)}}} \left(
\mathfrak{P}
\right)
\right)
\right]_{\varphi(s)}
&=
\left[
f_{s}^{(1)\flat} \left(
\mathrm{ip}_{\boldsymbol{\mathcal{A}}^{(1)}, s}^{(1,0)\sharp} \left(
\mathrm{sc}_{\boldsymbol{\mathcal{A}}^{(1)}, s}^{(0,1)} \left(
\mathfrak{P}
\right)
\right)
\right)
\right]_{\varphi(s)}
\tag{1}
\\
&=
\left[
\mathrm{ip}_{\boldsymbol{\mathcal{B}}^{(1)}, \varphi(s)}^{(1,0)\sharp} \left(
f_{s}^{(0)\sharp} \left(
\mathrm{sc}_{\boldsymbol{\mathcal{A}}^{(1)}, s}^{(0,1)} \left(
\mathfrak{P}
\right)
\right)
\right)
\right]_{\varphi(s)}
\tag{2}
\\
&=
\left[
\mathrm{ip}_{\boldsymbol{\mathcal{B}}^{(1)}, \varphi(s)}^{(1,0)\sharp} \left(
\mathrm{sc}_{\boldsymbol{\mathcal{B}}^{(1)}, \varphi(s)}^{(0,1)} \left(
f_{s}^{(1)\flat} \left(
\mathfrak{P}
\right)
\right)
\right)
\right]_{\varphi(s)}
\tag{3}
\\
&=
\left[
\mathrm{sc}_{\varphi(s)}^{0\mathbf{Pth}_{\boldsymbol{\mathcal{B}}^{(1)}}} \left(
f_{s}^{(1)\flat} \left(
\mathfrak{P}
\right)
\right)
\right]_{\varphi(s)}
\tag{4}
\\
&=
\left[
\mathrm{sc}_{s}^{0\mathbf{Pth}_{\boldsymbol{\mathcal{B}}^{(1)}}^{\mathbf{f}^{(1)}}} \left(
f_{s}^{(1)\flat} \left(
\mathfrak{P}
\right)
\right)
\right]_{\varphi(s)}
\tag{5}
\\
&=
\mathrm{sc}_{s}^{0[\mathbf{Pth}_{\boldsymbol{\mathcal{B}}^{(1)}}^{\mathbf{f}^{(1)}}]} \left(
\left[
f_{s}^{(1)\flat} \left(
\mathfrak{P}
\right)
\right]_{\varphi(s)}
\right).
\tag{6}
\end{align*}

The first equality unravels the interpretation of the operation symbol $\mathrm{sc}_{s}^{0}$ in the partial $\Sigma^{\boldsymbol{\mathcal{A}}^{(1)}}$-algebra $\mathbf{Pth}_{\boldsymbol{\mathcal{A}}^{(1)}}$, introduced in Proposition~\ref{PDPthDCatAlg}; 
the second equality follows from the definition of $f^{(1)\flat}$, introduced in Proposition~\ref{PPthExt}; 
the third equality follows from Proposition~\ref{PPthExt}; 
the fourth equality recovers the interpretation of the operation symbol $\mathrm{sc}_{\varphi(s)}^{0}$ in the partial $\Lambda^{\boldsymbol{\mathcal{B}}^{(1)}}$-algebra $\mathbf{Pth}_{\boldsymbol{\mathcal{B}}^{(1)}}$, according to Proposition~\ref{PDPthDCatAlg}; 
the fifth equality recovers the interpretation of the operation symbol $\mathrm{sc}_{s}^{0}$ in the partial $\Sigma^{\boldsymbol{\mathcal{A}}^{(1)}}$-algebra $\mathbf{Pth}_{\boldsymbol{\mathcal{B}}^{(1)}}^{\mathbf{f}^{(1)}}$, introduced in Proposition~\ref{PPthBCatAlg}; 
finally, the last equality follows from the fact that, according to Proposition~\ref{PKerCHBCatHom}, $\mathrm{pr}_{\boldsymbol{\mathcal{B}}^{(1)}, \varphi}^{\mathrm{Ker}(\mathrm{CH}^{(1)})}$ is a $\Sigma^{\boldsymbol{\mathcal{A}}^{(1)}}$-homomorphism from $\mathbf{Pth}_{\boldsymbol{\mathcal{B}}^{(1)}}^{\mathbf{f}^{(1)}}$ to $[\mathbf{Pth}_{\boldsymbol{\mathcal{B}}^{(1)}}^{\mathbf{f}^{(1)}}]$.

Hence, $\mathrm{pr}_{\boldsymbol{\mathcal{B}}^{(1)}, \varphi}^{\mathrm{Ker}(\mathrm{CH}^{(1)})} \circ f^{(1)\flat}$ is compatible with the $0$-source operation.

{\sffamily The mapping $\mathrm{pr}_{\boldsymbol{\mathcal{B}}^{(1)}, \varphi}^{\mathrm{Ker}(\mathrm{CH}^{(1)})} \circ f^{(1)\flat}$ is compatible with the $0$-target.}

Let $s$ be a sort in $S$ and let us consider the $0$-target operation symbol $\mathrm{tg}_{s}^{0}$ in $\Sigma^{\boldsymbol{\mathcal{A}}^{(1)}}_{s,s}$. Let $\mathfrak{P}$ be a path in $\mathrm{Pth}_{\boldsymbol{\mathcal{B}}^{(1)}, \varphi(s)}$, then the following equality holds
$$
\left[
f_{s}^{(1)\flat} \left(
\mathrm{tg}_{s}^{0\mathbf{Pth}_{\boldsymbol{\mathcal{A}}^{(1)}}} \left(
\mathfrak{P}
\right)
\right)
\right]_{\varphi(s)}
=
\mathrm{tg}_{s}^{0[\mathbf{Pth}_{\boldsymbol{\mathcal{B}}^{(1)}}^{\mathbf{f}^{(1)}}]} \left(
\left[
f_{s}^{(1)\flat} \left(
\mathfrak{P}
\right)
\right]_{\varphi(s)}
\right).
$$

The proof of this case is identical to that of the $0$-source.

Hence, $\mathrm{pr}_{\boldsymbol{\mathcal{B}}^{(1)}, \varphi}^{\mathrm{Ker}(\mathrm{CH}^{(1)})} \circ f^{(1)\flat}$ is compatible with the $0$-target operation.

{\sffamily The mapping $\mathrm{pr}_{\boldsymbol{\mathcal{B}}^{(1)}, \varphi}^{\mathrm{Ker}(\mathrm{CH}^{(1)})} \circ f^{(1)\flat}$ is compatible with the $0$-composition.}

Let $s$ be a sort in $S$ and let us consider the $0$-composition operation symbol $\circ_{s}^{0}$ in $\Sigma_{ss,s}^{\boldsymbol{\mathcal{A}}^{(1)}}$. Let $\mathfrak{P}$ and $\mathfrak{Q}$ be two paths in $\mathrm{Pth}_{\boldsymbol{\mathcal{A}}^{(1)}, s}$ such that
$$
\mathrm{sc}_{\boldsymbol{\mathcal{A}}^{(1)}, s}^{(0,1)}\left(\mathfrak{Q}\right)
=
\mathrm{tg}_{\boldsymbol{\mathcal{A}}^{(1)}, s}^{(0,1)}\left(\mathfrak{P}\right).
$$

We will prove that
$$
\left[
f_{s}^{(1)\flat}\left(
\mathfrak{Q}
\circ_{s}^{0\mathbf{Pth}_{\boldsymbol{\mathcal{A}}^{(1)}}}
\mathfrak{P}
\right)
\right]_{\varphi(s)}
=
\left[
f_{s}^{(1)\flat}\left(
\mathfrak{Q}
\right)
\right]_{\varphi(s)}
\circ_{s}^{0[\mathbf{Pth}_{\boldsymbol{\mathcal{B}}^{(1)}}^{\mathbf{f}^{(1)}}]}
\left[
f_{s}^{(1)\flat}\left(
\mathfrak{P}
\right)
\right]_{\varphi(s)}
$$
by Artinian induction of $(\coprod \mathrm{Pth}_{\boldsymbol{\mathcal{A}}^{(1)}}, \leq_{\mathrm{Pth}_{\boldsymbol{\mathcal{A}}^{(1)}}})$.

{\sf Base step of the Artinian induction.}

Let $(\mathfrak{Q}\circ_{s}^{0\mathbf{Pth}_{\boldsymbol{\mathcal{A}}^{(1)}}}\mathfrak{P}, s)$ be a minimal element in $(\coprod \mathrm{Pth}_{\boldsymbol{\mathcal{A}}^{(1)}}, \leq_{\mathrm{Pth}_{\boldsymbol{\mathcal{A}}^{(1)}}})$. Then, by Proposition~\ref{PMinimal}, the second order path $\mathfrak{Q} \circ_{s}^{0\mathbf{Pth}_{\boldsymbol{\mathcal{A}}^{(1)}}}\mathfrak{P}$ is either a $(1,0)$-identity path or an echelon. In any case, either $\mathfrak{Q}$ or $\mathfrak{P}$ must be a $(1,0)$-identity path. 

Assume that $\mathfrak{P}$ is a $(1,0)$-identity path.

Since $\mathfrak{Q}\circ_{s}^{0\mathbf{Pth}_{\boldsymbol{\mathcal{A}}^{(1)}}}\mathfrak{P}$ is defined, we have that $\mathrm{sc}_{\boldsymbol{\mathcal{A}}^{(1)}, s}^{(0,1)}(\mathfrak{Q}) = \mathrm{tg}_{\boldsymbol{\mathcal{A}}^{(1)}, s}^{(0,1)}(\mathfrak{P})$. Hence, $\mathfrak{P}$ is the $(1,0)$-identity path on the $(0,1)$-source of $\mathfrak{Q}$, i.e., $\mathfrak{P} = \mathrm{ip}_{\boldsymbol{\mathcal{A}}^{(1)}, s}^{(1,0)\sharp}(\mathrm{sc}_{s}^{(0,1)}(\mathfrak{P}))$. Hence, the $0$-composition $\mathfrak{Q} \circ_{s}^{0\mathbf{Pth}_{\boldsymbol{\mathcal{A}}^{(1)}}}\mathfrak{P}$ reduces to $\mathfrak{Q}$.

On the other hand, the following chain of equalities holds
\begin{align}
\left[
f_{s}^{(1)\flat}\left(
\mathfrak{P}
\right)
\right]_{\varphi(s)}
&=
\left[
f_{s}^{(1)\flat}\left(
\mathrm{ip}^{(1,0)\sharp}_{\boldsymbol{\mathcal{A}}^{(1)},s}\left(
\mathrm{sc}^{(0,1)}_{\boldsymbol{\mathcal{A}}^{(1)},s}\left(
\mathfrak{Q}
\right)
\right)
\right)
\right]_{\varphi(s)}
\tag{1}
\\
&=
\left[
\mathrm{ip}^{(1,0)\sharp}_{\boldsymbol{\mathcal{B}}^{(1)},\varphi(s)}\left(
f_{s}^{(1)@}\left(
\mathrm{sc}^{(0,1)}_{\boldsymbol{\mathcal{A}}^{(1)},s}\left(
\mathfrak{Q}
\right)
\right)
\right)
\right]_{\varphi(s)}
\tag{2}
\\
&=
\left[
\mathrm{ip}^{(1,0)\sharp}_{\boldsymbol{\mathcal{B}}^{(1)},\varphi(s)}\left(
\mathrm{sc}^{(0,1)}_{\boldsymbol{\mathcal{B}}^{(1)},\varphi(s)}\left(
f_{s}^{(1)\sharp}\left(
\mathfrak{Q}
\right)
\right)
\right)
\right]_{\varphi(s)}.
\tag{3}
\end{align}

The first equality follows from the fact that $\mathfrak{P} = \mathrm{ip}^{(1,0)\sharp}_{\boldsymbol{\mathcal{A}}^{(1)}, s}(\mathrm{sc}^{(0,1)}_{\boldsymbol{\mathcal{A}}^{(1)}, s}(\mathfrak{Q}))$;
the second equality follows from Proposition~\ref{PPthExt};
finally, the last equality follows from Proposition~\ref{PPthExt}.

Hence, the $0$-composition 
$
[
f_{s}^{(1)\flat}(
\mathfrak{Q}
)
]_{\varphi(s)}
\circ_{s}^{0[\mathbf{Pth}_{\boldsymbol{\mathcal{B}}^{(1)}}^{\mathbf{f}^{(1)}}]}
[
f_{s}^{(1)\flat}(
\mathfrak{P}
)
]_{\varphi(s)}
$
reduces to
$
[
f_{s}^{(1)\flat}(
\mathfrak{Q}
)
]_{\varphi(s)}
$.

All in all, we conclude that
\allowdisplaybreaks
\begin{align}
\left[
f_{s}^{(1)\flat}\left(
\mathfrak{Q}
\circ_{s}^{0\mathbf{Pth}_{\boldsymbol{\mathcal{A}}^{(1)}}}
\mathfrak{P}
\right)
\right]_{\varphi(s)}
&=
\left[
f_{s}^{(1)\flat}\left(
\mathfrak{Q}
\right)
\right]_{\varphi(s)}
\tag{1}
\\
&=
\left[
f_{s}^{(1)\flat}\left(
\mathfrak{Q}
\right)
\right]_{\varphi(s)}
\circ_{s}^{0[\mathbf{Pth}_{\boldsymbol{\mathcal{B}}^{(1)}}^{\mathbf{f}^{(1)}}]}
\left[
f_{s}^{(1)\flat}\left(
\mathfrak{P}
\right)
\right]_{\varphi(s)}
\tag{2}
\end{align}

The first equality follows from the fact that the $0$-composition $\mathfrak{Q} \circ_{s}^{0\mathbf{Pth}_{\boldsymbol{\mathcal{A}}^{(1)}}}\mathfrak{P}$ reduces to $\mathfrak{Q}$; 
the second equality follows from the fact that the $0$-composition 
$
[
f_{s}^{(1)\flat}(
\mathfrak{Q}
)
]_{\varphi(s)}
\circ_{s}^{0[\mathbf{Pth}_{\boldsymbol{\mathcal{B}}^{(1)}}^{\mathbf{f}^{(1)}}]}
[
f_{s}^{(1)\flat}(
\mathfrak{P}
)
]_{\varphi(s)}
$
reduces to
$
[
f_{s}^{(1)\flat}(
\mathfrak{Q}
)
]_{\varphi(s)}
$.

The case where $\mathfrak{Q}$ is a $(1,0)$-identity path follows similarly.

This concludes the base step of the Artinian induction.

{\sf Inductive step of the Artinian induction.}

Let $(\mathfrak{Q} \circ_{s}^{0\mathbf{Pth}_{\boldsymbol{\mathcal{A}}^{(1)}}} \mathfrak{P}, s)$ be a non-minimal element in $(\coprod\mathrm{Pth}_{\boldsymbol{\mathcal{A}}^{(1)}}, \leq_{\mathbf{Pth}_{\boldsymbol{\mathcal{A}}^{(1)}}})$. Let us suppose that, for every sort $t \in S$ and every path $\mathfrak{Q}' \circ_{t}^{0\mathbf{Pth}_{\boldsymbol{\mathcal{A}}^{(1)}}} \mathfrak{P}' \in \mathrm{Pth}_{\boldsymbol{\mathcal{A}}^{(1)}, t}$, if $(\mathfrak{Q}' \circ_{t}^{0\mathbf{Pth}_{\boldsymbol{\mathcal{A}}^{(1)}}} \mathfrak{P}', t) <_{\mathbf{Pth}_{\boldsymbol{\mathcal{A}}^{(1)}}} (\mathfrak{Q} \circ_{s}^{0\mathbf{Pth}_{\boldsymbol{\mathcal{A}}^{(1)}}} \mathfrak{P}, s)$, then
$$
\left[
f_{t}^{(1)\flat}\left(
\mathfrak{Q}'
\circ_{t}^{0\mathbf{Pth}_{\boldsymbol{\mathcal{A}}^{(1)}}}
\mathfrak{P}'
\right)
\right]_{\varphi(t)}
=
\left[
f_{t}^{(1)\flat}\left(
\mathfrak{Q}'
\right)
\right]_{\varphi(t)}
\circ_{t}^{0[\mathbf{Pth}_{\boldsymbol{\mathcal{B}}^{(1)}}^{\mathbf{f}^{(1)}}]} 
\left[
f_{t}^{(1)\flat}\left(
\mathfrak{P}'
\right)
\right]_{\varphi(t)}.
$$

Since $\mathfrak{Q} \circ_{s}^{0\mathbf{Pth}_{\boldsymbol{\mathcal{A}}^{(1)}}} \mathfrak{P}$ is a non-minimal element in $(\coprod\mathrm{Pth}_{\boldsymbol{\mathcal{A}}^{(1)}}, \leq_{\mathbf{Pth}_{\boldsymbol{\mathcal{A}}^{(1)}}})$ and we can assume that neither $\mathfrak{Q}$ nor $\mathfrak{P}$ are $(1,0)$-identity paths because this case already has been considered, we have, by Lemma~\ref{LOrdI}, that $\mathfrak{Q} \circ_{s}^{0\mathbf{Pth}_{\boldsymbol{\mathcal{A}}^{(1)}}} \mathfrak{P}$ is either~(1) a path of length strictly grater than one containing at least one echelon or~(2) an echelonless path.

If~(1), then let $i\in\bb{\mathfrak{Q} \circ_{s}^{0\mathbf{Pth}_{\boldsymbol{\mathcal{A}}^{(1)}}} \mathfrak{P}}$ be the first index for which the one-step subpath $(\mathfrak{Q} \circ_{s}^{0\mathbf{Pth}_{\boldsymbol{\mathcal{A}}^{(1)}}} \mathfrak{P})^{i,i}$ of $\mathfrak{Q} \circ_{s}^{0\mathbf{Pth}_{\boldsymbol{\mathcal{A}}^{(1)}}} \mathfrak{P}$ is an echelon. We distinguish the cases~(1.1) $i=0$ and~(1.2) $i>0$.

If~(1.1), i.e., $i=0$, since we are assuming that $\mathfrak{P}$ is not a (1,0)-identity path, we have that $\mathfrak{P}$ has an echelon on its first step. Then it could be the case that either~(1.1.1) $\mathfrak{P}$ is an echelon or~(1.1.2) $\mathfrak{P}$ is a path of length strictly greater than one containing an echelon on its first step.

If~(1.1.1), then, regarding the paths $\mathfrak{Q}$ and $\mathfrak{P}$, we have that
\begin{enumerate}
\item[(i)]
$\mathfrak{P}$ is an echelon.
\item[(ii)]
$\mathfrak{Q}$ is a path.
\end{enumerate}

So, considering the foregoing, we can affirm that
\allowdisplaybreaks
\begin{align*}
\left[
f_{s}^{(1)\flat}\left(
\mathfrak{Q}
\circ_{s}^{0\mathbf{Pth}_{\boldsymbol{\mathcal{A}}^{(1)}}}
\mathfrak{P}
\right)
\right]_{\varphi(s)}
&=
\left[
f_{s}^{(1)\flat}\left(
\mathfrak{Q}
\right)
\circ_{\varphi(s)}^{0\mathbf{Pth}_{\boldsymbol{\mathcal{B}}^{(1)}}}
f_{s}^{(1)\flat}\left(
\mathfrak{P}
\right)
\right]_{\varphi(s)}
\tag{1}
\\
&=
\left[
f_{s}^{(1)\flat}\left(
\mathfrak{Q}
\right)
\circ_{s}^{0\mathbf{Pth}_{\boldsymbol{\mathcal{B}}^{(1)}}^{\mathbf{f}^{(1)}}}
f_{s}^{(1)\flat}\left(
\mathfrak{P}
\right)
\right]_{\varphi(s)}
\tag{2}
\\
&=
\left[
f_{s}^{(1)\flat}\left(
\mathfrak{Q}
\right)
\right]_{\varphi(s)}
\circ_{s}^{0[\mathbf{Pth}_{\boldsymbol{\mathcal{B}}^{(1)}}^{\mathbf{f}^{(1)}}]}
\left[
f_{s}^{(1)\flat}\left(
\mathfrak{P}
\right)
\right]_{\varphi(s)}.
\tag{3}
\end{align*}

The first equality unravels the definition of the path extension mapping $f^{(1)\flat}$ introduced in Proposition~\ref{PPthExt};
the second equality recovers the interpretation of the operation symbol $\circ_{s}^{0}$ in the partial $\Sigma^{\boldsymbol{\mathcal{A}}^{(1)}}$-algebra $\mathbf{Pth}_{\boldsymbol{\mathcal{B}}^{(1)}}^{\mathbf{f}^{(1)}}$ introduced in Proposition~\ref{PPthBCatAlg};
finally, the last equality recovers the interpretation of the operation symbol $\circ_{s}^{0}$ in the partial $\Sigma^{\boldsymbol{\mathcal{A}}^{(1)}}$-algebra $[\mathbf{Pth}_{\boldsymbol{\mathcal{B}}^{(1)}}^{\mathbf{f}^{(1)}}]$ introduced in Proposition~\ref{PQPthBCatAlg}.

If~(1.1.2), then, regarding the paths $\mathfrak{Q}$ and $\mathfrak{P}$, we have that
\begin{enumerate}
\item[(i)]
$\mathfrak{P}$ is a path of length strictly grater than one containing an echelon on its first step.
\item[(ii)]
$\mathfrak{Q}$ is a path.
\end{enumerate}

From (i) and taking into account the definition of $f^{(1)\flat}$ introduced in Proposition~\ref{PPthExt}, we have that the value of the path extension mapping $f^{(1)\flat}$ at $\mathfrak{P}$ is given by
$$
f_{s}^{(1)\flat}\left(
\mathfrak{P}
\right)
=
f_{s}^{(1)\flat}\left(
\mathfrak{P}^{1,\bb{\mathfrak{P}}-1}
\right)
\circ^{1\mathbf{Pth}_{\boldsymbol{\mathcal{B}}^{(1)}}}_{\varphi(s)}
f_{s}^{(1)\flat}\left(
\mathfrak{P}^{0,0}
\right).
$$

Since $(\mathfrak{Q} \circ_{s}^{0\mathbf{Pth}_{\boldsymbol{\mathcal{A}}^{(1)}}} \mathfrak{P}^{1,\bb{\mathfrak{P}}-1}, s)$$\prec_{\mathbf{Pth}_{\boldsymbol{\mathcal{A}}^{(1)}}}$$(\mathfrak{Q} \circ_{s}^{0\mathbf{Pth}_{\boldsymbol{\mathcal{A}}^{(1)}}} \mathfrak{P}, s)$ we have, by induction, that
\begin{align*}
&\left[
f_{s}^{(1)\flat}\left(
\mathfrak{Q}
\circ_{s}^{0\mathbf{Pth}_{\boldsymbol{\mathcal{A}}^{(1)}}}
\mathfrak{P}^{1,\bb{\mathfrak{P}}-1}
\right)
\right]_{\varphi(s)}
=
\left[
f_{s}^{(1)\flat}\left(
\mathfrak{Q}
\right)
\right]_{\varphi(s)}
\circ_{s}^{0[\mathbf{Pth}_{\boldsymbol{\mathcal{B}}^{(1)}}^{\mathbf{f}^{(1)}}]}
\left[
f_{s}^{(1)\flat}\left(
\mathfrak{P}^{1,\bb{\mathfrak{P}}-1}
\right)
\right]_{\varphi(s)}.
\end{align*}

So, considering the foregoing, we can affirm that
\begin{flushleft}
$
\left[
f_{s}^{(1)\flat}\left(
\mathfrak{Q} \circ_{s}^{0\mathbf{Pth}_{\boldsymbol{\mathcal{A}}^{(1)}}} \mathfrak{P}
\right)
\right]_{\varphi(s)}
$
\allowdisplaybreaks
\begin{align*}
&=
\left[
f_{s}^{(1)\flat}\left(
\mathfrak{Q} \circ_{s}^{0\mathbf{Pth}_{\boldsymbol{\mathcal{A}}^{(1)}}} \mathfrak{P}^{1,\bb{\mathfrak{P}}-1}
\right)
\circ_{\varphi(s)}^{0\mathbf{Pth}_{\boldsymbol{\mathcal{B}}^{(1)}}}
f_{s}^{(1)\flat}\left(
\mathfrak{P}^{0,0}
\right)
\right]_{\varphi(s)}
\tag{1}
\\
&=
\left[
f_{s}^{(1)\flat}\left(
\mathfrak{Q} \circ_{s}^{0\mathbf{Pth}_{\boldsymbol{\mathcal{A}}^{(1)}}} \mathfrak{P}^{1,\bb{\mathfrak{P}}-1}
\right)
\circ_{s}^{0\mathbf{Pth}_{\boldsymbol{\mathcal{B}}^{(1)}}^{\mathbf{f}^{(1)}}}
f_{s}^{(1)\flat}\left(
\mathfrak{P}^{0,0}
\right)
\right]_{\varphi(s)}
\tag{2}
\\
&=
\left[
f_{s}^{(1)\flat}\left(
\mathfrak{Q} \circ_{s}^{1\mathbf{Pth}_{\boldsymbol{\mathcal{A}}^{(1)}}} \mathfrak{P}^{1,\bb{\mathfrak{P}}-1}
\right)
\right]_{\varphi(s)}
\circ_{s}^{0[\mathbf{Pth}_{\boldsymbol{\mathcal{B}}^{(1)}}^{\mathbf{f}^{(1)}}]}
\left[
f_{s}^{(1)\flat}\left(
\mathfrak{P}^{0,0}
\right)
\right]_{\varphi(s)}
\tag{3}
\\
&=
\left(
\left[
f_{s}^{(1)\flat}\left(
\mathfrak{Q}
\right)
\right]_{\varphi(s)}
\circ_{s}^{0[\mathbf{Pth}_{\boldsymbol{\mathcal{B}}^{(1)}}^{\mathbf{f}^{(1)}}]} 
\left[
f_{s}^{(1)\flat}\left(
\mathfrak{P}^{1,\bb{\mathfrak{P}}-1}
\right)
\right]_{\varphi(s)}
\right)
\\
&\hspace{6.7cm}
\circ_{s}^{0[\mathbf{Pth}_{\boldsymbol{\mathcal{B}}^{(1)}}^{\mathbf{f}^{(1)}}]}
\left[
f_{s}^{(1)\flat}\left(
\mathfrak{P}^{0,0}
\right)
\right]_{\varphi(s)}
\tag{4}
\\
&=
\left[
f_{s}^{(1)\flat}\left(
\mathfrak{Q}
\right)
\right]_{\varphi(s)}
\circ_{s}^{0[\mathbf{Pth}_{\boldsymbol{\mathcal{B}}^{(1)}}^{\mathbf{f}^{(1)}}]} 
\\
&\hspace{3cm}
\left(
\left[
f_{s}^{(1)\flat}\left(
\mathfrak{P}^{1,\bb{\mathfrak{P}}-1}
\right)
\right]_{\varphi(s)}
\circ_{s}^{0[\mathbf{Pth}_{\boldsymbol{\mathcal{B}}^{(1)}}^{\mathbf{f}^{(1)}}]}
\left[
f_{s}^{(1)\flat}\left(
\mathfrak{P}^{0,0}
\right)
\right]_{\varphi(s)}
\right)
\tag{5}
\\
&=
\left[
f_{s}^{(1)\flat}\left(
\mathfrak{Q}
\right)
\right]_{\varphi(s)}
\circ_{s}^{0[\mathbf{Pth}_{\boldsymbol{\mathcal{B}}^{(1)}}^{\mathbf{f}^{(1)}}]}
\left[
f_{s}^{(1)\flat}\left(
\mathfrak{P}^{1,\bb{\mathfrak{P}}-1}
\right)
\circ_{s}^{0\mathbf{Pth}_{\boldsymbol{\mathcal{B}}^{(1)}}^{\mathbf{f}^{(1)}}}
f_{s}^{(1)\flat}\left(
\mathfrak{P}^{0,0}
\right)
\right]_{\varphi(s)}
\tag{6}
\\
&=
\left[
f_{s}^{(1)\flat}\left(
\mathfrak{Q}
\right)
\right]_{\varphi(s)}
\circ_{s}^{0[\mathbf{Pth}_{\boldsymbol{\mathcal{B}}^{(1)}}^{\mathbf{f}^{(1)}}]} 
\left[
f_{s}^{(1)\flat}\left(
\mathfrak{P}^{1,\bb{\mathfrak{P}}-1}
\right)
\circ_{\varphi(s)}^{0\mathbf{Pth}_{\boldsymbol{\mathcal{B}}^{(1)}}}
f_{s}^{(1)\flat}\left(
\mathfrak{P}^{0,0}
\right)
\right]_{\varphi(s)}
\tag{7}
\\
&=
\left[
f_{s}^{(1)\flat}\left(
\mathfrak{Q}
\right)
\right]_{\varphi(s)}
\circ_{s}^{0[\mathbf{Pth}_{\boldsymbol{\mathcal{B}}^{(1)}}^{\mathbf{f}^{(1)}}]} 
\left[
f_{s}^{(1)\flat}\left(
\mathfrak{P}^{1,\bb{\mathfrak{P}}-1}
\circ_{s}^{0\mathbf{Pth}_{\boldsymbol{\mathcal{A}}^{(1)}}}
\mathfrak{P}^{0,0}
\right)
\right]_{\varphi(s)}
\tag{8}
\\
&=
\left[
f_{s}^{(1)\flat}\left(
\mathfrak{Q}
\right)
\right]_{\varphi(s)}
\circ_{s}^{0[\mathbf{Pth}_{\boldsymbol{\mathcal{B}}^{(1)}}^{\mathbf{f}^{(1)}}]} 
\left[
f_{s}^{(1)\flat}\left(
\mathfrak{P}
\right)
\right]_{\varphi(s)}.
\tag{9}
\end{align*}
\end{flushleft}

The first equality unravels the definition of the path extension mapping $f^{(1)\flat}$ introduced in Proposition~\ref{PPthExt};
the second equality recovers the definition of the operation symbol $\circ_{s}^{0}$ in the partial $\Sigma^{\boldsymbol{\mathcal{A}}^{(1)}}$-algebra $\mathbf{Pth}_{\boldsymbol{\mathcal{B}}^{(1)}}^{\mathbf{f}^{(1)}}$ introduced in Proposition~\ref{PPthBCatAlg};
the third equality recovers the definition of the operation symbol $\circ_{s}^{0}$ in the partial $\Sigma^{\boldsymbol{\mathcal{A}}^{(1)}}$-algebra $[\mathbf{Pth}_{\boldsymbol{\mathcal{B}}^{(1)}}^{\mathbf{f}^{(1)}}]$ introduced in Proposition~\ref{PQPthBCatAlg};
the fourth equality follows by Artinian induction;
the fifth equality follows from the fact that, according to Proposition~\ref{PPthComp}, the interpretation of the partial operation $\circ_{s}^{0}$ is associative;
the sixth equality unravels the definition of the operation symbol $\circ_{s}^{0}$ in the partial $\Sigma^{\boldsymbol{\mathcal{A}}^{(1)}}$-algebra $[\mathbf{Pth}_{\boldsymbol{\mathcal{B}}^{(1)}}^{\mathbf{f}^{(1)}}]$ introduced in Proposition~\ref{PQPthBCatAlg};
The seventh equality unravels the definition of the operation symbol $\circ_{s}^{0}$ in the partial $\Sigma^{\boldsymbol{\mathcal{A}}^{(1)}}$-algebra $\mathbf{Pth}_{\boldsymbol{\mathcal{B}}^{(1)}}^{\mathbf{f}^{(1)}}$ introduced in Proposition~\ref{PPthBCatAlg};
the eight equality recovers the definition of $f^{(1)\flat}$ since $\mathfrak{P}$ is a path of length strictly greater than one containing an echelon on its first step;
finally, the last equality recovers the definition of $\mathfrak{P}$.

This finishes the case $i = 0$.

For the case~$(1.2)$, i.e., if $i\neq 0$, since $\bb{\mathfrak{Q} \circ_{s}^{0\mathbf{Pth}_{\boldsymbol{\mathcal{A}}^{(1)}}} \mathfrak{P}} = \bb{\mathfrak{Q}} + \bb{\mathfrak{P}}$, then either~$(1.2.1)$ $i\in\bb{\mathfrak{P}}$ or~$(1.2.2)$ $i\in[\bb{P}, \bb{\mathfrak{Q} \circ_{s}^{0\mathbf{Pth}_{\boldsymbol{\mathcal{A}}^{(1)}}} \mathfrak{P}}-1]$.

If~(1.2.1), i.e., if we find ourselves in the case where $\mathfrak{Q} \circ_{s}^{0\mathbf{Pth}_{\boldsymbol{\mathcal{A}}^{(1)}}} \mathfrak{P}$ is a path that is not echelonless and $i > 0$, $i \in \bb{\mathfrak{P}}$ is the first index for which the one step subpath $(\mathfrak{Q} \circ_{s}^{0\mathbf{Pth}_{\boldsymbol{\mathcal{A}}^{(1)}}} \mathfrak{P})^{i,i}$ is an echelon then, regarding the paths $\mathfrak{Q}$ and $\mathfrak{P}$, we have that
\begin{enumerate}
\item[(i)]
$\mathfrak{P}$ is a path of length strictly grater than one containing an echelon on a step different from zero.
\item[(ii)]
$\mathfrak{Q}$ is a path.
\end{enumerate}

From (i) and taking into account the definition of $f^{(1)\flat}$ introduced in Proposition~\ref{PPthExt}, we have that the value of the path extension mapping $f^{(1)\flat}$ at $\mathfrak{P}$ is given by
$$
f_{s}^{(1)\flat}\left(
\mathfrak{P}
\right)
=
f_{s}^{(1)\flat}\left(
\mathfrak{P}^{i,\bb{\mathfrak{P}}-1}
\right)
\circ^{0\mathbf{Pth}_{\boldsymbol{\mathcal{B}}^{(1)}}}_{\varphi(s)}
f_{s}^{(1)\flat}\left(
\mathfrak{P}^{0,i-1}
\right).
$$

Since $(\mathfrak{Q} \circ_{s}^{0\mathbf{Pth}_{\boldsymbol{\mathcal{A}}^{(1)}}} \mathfrak{P}^{i,\bb{\mathfrak{P}}-1}, s)$$\prec_{\mathbf{Pth}_{\boldsymbol{\mathcal{A}}^{(1)}}}$$(\mathfrak{Q} \circ_{s}^{0\mathbf{Pth}_{\boldsymbol{\mathcal{A}}^{(1)}}} \mathfrak{P}, s)$ we have, by induction, that
\begin{align*}
&\left[
f_{s}^{(1)\flat}\left(
\mathfrak{Q}
\circ_{s}^{0\mathbf{Pth}_{\boldsymbol{\mathcal{A}}^{(1)}}}
\mathfrak{P}^{i,\bb{\mathfrak{P}}-1}
\right)
\right]_{\varphi(s)}
\\
&\hspace{3cm}
=
\left[
f_{s}^{(1)\flat}\left(
\mathfrak{Q}
\right)
\right]_{\varphi(s)}
\circ_{s}^{0[\mathbf{Pth}_{\boldsymbol{\mathcal{B}}^{(1)}}^{\mathbf{f}^{(1)}}]}
\left[
f_{s}^{(1)\flat}\left(
\mathfrak{P}^{i,\bb{\mathfrak{P}}-1}
\right)
\right]_{\varphi(s)}.
\end{align*}

So, considering the foregoing, we can affirm that
\begin{flushleft}
$
\left[
f_{s}^{(1)\flat}\left(
\mathfrak{Q} \circ_{s}^{0\mathbf{Pth}_{\boldsymbol{\mathcal{A}}^{(1)}}} \mathfrak{P}
\right)
\right]_{\varphi(s)}
$
\allowdisplaybreaks
\begin{align*}
&=
\left[
f_{s}^{(1)\flat}\left(
\mathfrak{Q} \circ_{s}^{0\mathbf{Pth}_{\boldsymbol{\mathcal{A}}^{(1)}}} \mathfrak{P}^{i,\bb{\mathfrak{P}}-1}
\right)
\circ_{\varphi(s)}^{0\mathbf{Pth}_{\boldsymbol{\mathcal{B}}^{(1)}}}
f_{s}^{(1)\flat}\left(
\mathfrak{P}^{0,i-1}
\right)
\right]_{\varphi(s)}
\tag{1}
\\
&=
\left[
f_{s}^{(1)\flat}\left(
\mathfrak{Q} \circ_{s}^{1\mathbf{Pth}_{\boldsymbol{\mathcal{A}}^{(1)}}} \mathfrak{P}^{(1)i,\bb{\mathfrak{P}}-1}
\right)
\circ_{s}^{0\mathbf{Pth}_{\boldsymbol{\mathcal{B}}^{(1)}}^{\mathbf{f}^{(1)}}}
f_{s}^{(1)\flat}\left(
\mathfrak{P}^{(1)0,i-1}
\right)
\right]_{\varphi(s)}
\tag{2}
\\
&=
\left[
f_{s}^{(1)\flat}\left(
\mathfrak{Q} \circ_{s}^{0\mathbf{Pth}_{\boldsymbol{\mathcal{A}}^{(1)}}} \mathfrak{P}^{(1)i,\bb{\mathfrak{P}}-1}
\right)
\right]_{\varphi(s)}
\circ_{s}^{1[\mathbf{Pth}_{\boldsymbol{\mathcal{B}}^{(1)}}^{\mathbf{f}^{(1)}}]}
\left[
f_{s}^{(1)\flat}\left(
\mathfrak{P}^{0,i-1}
\right)
\right]_{\varphi(s)}
\tag{3}
\\
&=
\left(
\left[
f_{s}^{(1)\flat}\left(
\mathfrak{Q}
\right)
\right]_{\varphi(s)}
\circ_{s}^{0[\mathbf{Pth}_{\boldsymbol{\mathcal{B}}^{(1)}}^{\mathbf{f}^{(1)}}]} 
\left[
f_{s}^{(1)\flat}\left(
\mathfrak{P}^{i,\bb{\mathfrak{P}}-1}
\right)
\right]_{\varphi(s)}
\right)
\\
&\hspace{6.8cm}
\circ_{s}^{0[\mathbf{Pth}_{\boldsymbol{\mathcal{B}}^{(1)}}^{\mathbf{f}^{(1)}}]}
\left[
f_{s}^{(1)\flat}\left(
\mathfrak{P}^{0,i-1}
\right)
\right]_{\varphi(s)}
\tag{4}
\\
&=
\left[
f_{s}^{(1)\flat}\left(
\mathfrak{Q}
\right)
\right]_{\varphi(s)}
\circ_{s}^{0[\mathbf{Pth}_{\boldsymbol{\mathcal{B}}^{(1)}}^{\mathbf{f}^{(1)}}]} 
\\
&\hspace{3cm}
\left(
\left[
f_{s}^{(1)\flat}\left(
\mathfrak{P}^{i,\bb{\mathfrak{P}}-1}
\right)
\right]_{\varphi(s)}
\circ_{s}^{0[\mathbf{Pth}_{\boldsymbol{\mathcal{B}}^{(1)}}^{\mathbf{f}^{(1)}}]}
\left[
f_{s}^{(1)\flat}\left(
\mathfrak{P}^{0,i-1}
\right)
\right]_{\varphi(s)}
\right)
\tag{5}
\\
&=
\left[
f_{s}^{(1)\flat}\left(
\mathfrak{Q}
\right)
\right]_{\varphi(s)}
\circ_{s}^{0[\mathbf{Pth}_{\boldsymbol{\mathcal{B}}^{(1)}}^{\mathbf{f}^{(1)}}]} 
\left[
f_{s}^{(1)\flat}\left(
\mathfrak{P}^{i,\bb{\mathfrak{P}}-1}
\right)
\circ_{s}^{1\mathbf{Pth}_{\boldsymbol{\mathcal{B}}^{(1)}}^{\mathbf{f}^{(1)}}}
f_{s}^{(1)\flat}\left(
\mathfrak{P}^{0,i-1}
\right)
\right]_{\varphi(s)}
\tag{6}
\\
&=
\left[
f_{s}^{(1)\flat}\left(
\mathfrak{Q}
\right)
\right]_{\varphi(s)}
\circ_{s}^{0[\mathbf{Pth}_{\boldsymbol{\mathcal{B}}^{(1)}}^{\mathbf{f}^{(1)}}]}
\left[
f_{s}^{(1)\flat}\left(
\mathfrak{P}^{i,\bb{\mathfrak{P}}-1}
\right)
\circ_{\varphi(s)}^{0\mathbf{Pth}_{\boldsymbol{\mathcal{B}}^{(1)}}}
f_{s}^{(1)\flat}\left(
\mathfrak{P}^{0,i-1}
\right)
\right]_{\varphi(s)}
\tag{7}
\\
&=
\left[
f_{s}^{(1)\flat}\left(
\mathfrak{Q}
\right)
\right]_{\varphi(s)}
\circ_{s}^{0[\mathbf{Pth}_{\boldsymbol{\mathcal{B}}^{(1)}}^{\mathbf{f}^{(1)}}]}
\left[
f_{s}^{(1)\flat}\left(
\mathfrak{P}^{i,\bb{\mathfrak{P}}-1}
\circ_{s}^{0\mathbf{Pth}_{\boldsymbol{\mathcal{A}}^{(1)}}}
\mathfrak{P}^{0,i-1}
\right)
\right]_{\varphi(s)}
\tag{8}
\\
&=
\left[
f_{s}^{(1)\flat}\left(
\mathfrak{Q}
\right)
\right]_{\varphi(s)}
\circ_{s}^{0[\mathbf{Pth}_{\boldsymbol{\mathcal{B}}^{(1)}}^{\mathbf{f}^{(1)}}]} 
\left[
f_{s}^{(1)\flat}\left(
\mathfrak{P}
\right)
\right]_{\varphi(s)}.
\tag{9}
\end{align*}
\end{flushleft}

The first equality unravels the definition of the path extension mapping $f^{(1)\flat}$ introduced in Proposition~\ref{PPthExt};
the second equality recovers the definition of the operation symbol $\circ_{s}^{0}$ in the partial $\Sigma^{\boldsymbol{\mathcal{A}}^{(1)}}$-algebra $\mathbf{Pth}_{\boldsymbol{\mathcal{B}}^{(1)}}^{\mathbf{f}^{(1)}}$ introduced in Proposition~\ref{PPthBCatAlg};
the third equality recovers the definition of the operation symbol $\circ_{s}^{0}$ in the partial $\Sigma^{\boldsymbol{\mathcal{A}}^{(1)}}$-algebra $[\mathbf{Pth}_{\boldsymbol{\mathcal{B}}^{(1)}}^{\mathbf{f}^{(1)}}]$ introduced in Proposition~\ref{PQPthBCatAlg};
the fourth equality follows by Artinian induction;
the fifth equality follows from the fact that, according to Proposition~\ref{PPthComp}, the interpretation of the partial operation $\circ_{s}^{0}$ is associative;
the sixth equality unravels the definition of the operation symbol $\circ_{s}^{0}$ in the partial $\Sigma^{\boldsymbol{\mathcal{A}}^{(1)}}$-algebra $[\mathbf{Pth}_{\boldsymbol{\mathcal{B}}^{(1)}}^{\mathbf{f}^{(1)}}]$ introduced in Proposition~\ref{PQPthBCatAlg};
The seventh equality unravels the definition of the operation symbol $\circ_{s}^{0}$ in the partial $\Sigma^{\boldsymbol{\mathcal{A}}^{(1)}}$-algebra $\mathbf{Pth}_{\boldsymbol{\mathcal{B}}^{(1)}}^{\mathbf{f}^{(1)}}$ introduced in Proposition~\ref{PPthBCatAlg};
the eight equality recovers the definition of $f^{(1)\flat}$ since $\mathfrak{P}$ is a path of length strictly greater than one containing an echelon on its first step;
finally, the last equality recovers the definition of $\mathfrak{P}$.

If~(1.2.2), i.e, $i \neq 0$ and $i \in [\bb{\mathfrak{P}}, \bb{\mathfrak{Q} \circ_{s}^{0\mathbf{Pth}_{\boldsymbol{\mathcal{A}}^{(1)}} \mathfrak{P}}-1}]$, then $\mathfrak{Q}$ is not an $(1,0)$-identity path containing an echelon, whilst $\mathfrak{P}$ is an echelonless path.

We will distinguish three cases according to whether (1.2.2.1) $\mathfrak{Q}$ is an echelon; (1.2.2.2) $\mathfrak{Q}$ is a path of length strictly greater than one containing an echelon of its first step or (1.2.2.3) $\mathfrak{Q}$ is a path of length strictly greater than one containing an echelon on a step different from zero. These cases can be proved using a similar argument to those three cases presented above. We leave the details for the interested reader.

This finishes the case $i > 0$.

This completes case $(1)$.

If~(2), i.e., if $\mathfrak{Q} \circ_{s}^{0\mathbf{Pth}_{\boldsymbol{\mathcal{A}}^{(1)}}} \mathfrak{P}$ is an echelonless path then, regarding the paths $\mathfrak{Q}$ and $\mathfrak{P}$, we have that
\begin{enumerate}
\item[(i)]
$\mathfrak{P}$ is an echelonless path.
\item[(ii)]
$\mathfrak{Q}$ is an echelonless path.
\end{enumerate}

Since (i), we have by Lemma~\ref{LPthHeadCt}, that there exists a unique word $\mathbf{s} \in S^{\star} - \{\lambda\}$ and a unique operation symbol $\sigma \in \Sigma_{\mathbf{s}, s}$ associated to $\mathfrak{P}$. Let $(\mathfrak{P}_{j})_{j \in \bb{s}}$ be the family of paths we can extract from $\mathfrak{P}$ in virtue of Lemma~\ref{LPthExtract}. Then according to Proposition~\ref{PPthExt}, we have that the value of the path extension mapping $f^{(1)\flat}$ at $\mathfrak{P}$ is given by
$$
f^{(1)\flat}_{s}(\mathfrak{P})
=
\sigma^{\mathbf{Pth}_{\boldsymbol{\mathcal{B}}^{(1)}}^{\mathbf{f}^{(1)}(0,1)}} \left(
\left( 
f^{(1)\flat}_{s_{j}} \left(
\mathfrak{P}_{j}
\right)
\right)_{j \in \bb{\mathbf{s}}} 
\right).
$$

Since (ii), we have by Lemma~\ref{LPthHeadCt}, that there exists a unique word $\mathbf{s} \in S^{\star} - \{\lambda\}$ and a unique operation symbol $\sigma \in \Sigma_{\mathbf{s}, s}$ associated to $\mathfrak{Q}$. Let $(\mathfrak{Q}_{j})_{j \in \bb{s}}$ be the family of paths we can extract from $\mathfrak{Q}$ in virtue of Lemma~\ref{LPthExtract}. Then according to Proposition~\ref{PPthExt}, we have that the value of the path extension mapping $f^{(1)\flat}$ at $\mathfrak{Q}$ is given by
$$
f^{(1)\flat}_{s}(\mathfrak{Q})
=
\sigma^{\mathbf{Pth}_{\boldsymbol{\mathcal{B}}^{(1)}}^{\mathbf{f}^{(1)}(0,1)}} \left(
\left( 
f^{(1)\flat}_{s_{j}} \left(
\mathfrak{Q}_{j}
\right)
\right)_{j \in \bb{\mathbf{s}}} 
\right).
$$

Note that the operation symbol $\sigma$ is the same as in case (i), since $\mathfrak{Q}\circ_{s}^{0\mathbf{Pth}_{\boldsymbol{\mathcal{A}}^{(1)}}}\mathfrak{P}$ is an echelonless path by hypothesis and, thus, head-constant by Lemma~\ref{LPthHeadCt}.

Let us consider $((\mathfrak{Q} \circ_{s}^{0\mathbf{Pth}_{\boldsymbol{\mathcal{A}}^{(1)}}} \mathfrak{P})_{j})_{j\in\bb{\mathbf{s}}}$ the family of paths we can extract, in virtue of Lemma~\ref{LPthExtract}, from $\mathfrak{Q} \circ_{s}^{0\mathbf{Pth}_{\boldsymbol{\mathcal{A}}^{(1)}}} \mathfrak{P}$. Let us note that, for every $j\in\bb{\mathbf{s}}$, it is the case that
$$
(\mathfrak{Q} \circ_{s}^{0\mathbf{Pth}_{\boldsymbol{\mathcal{A}}^{(1)}}} \mathfrak{P})_{j}
=
\mathfrak{Q}_{j} \circ_{s_{j}}^{0\mathbf{Pth}_{\boldsymbol{\mathcal{A}}^{(1)}}} \mathfrak{P}_{j}.
$$
Then according to Proposition~\ref{PPthExt}, we have that the value of the path extension mapping of $\mathbf{f}^{(1)}$ at $\mathfrak{Q} \circ_{s}^{0\mathbf{Pth}_{\boldsymbol{\mathcal{A}}^{(1)}}} \mathfrak{P}$ is given by
$$
f_{s}^{(1)\flat} \left(
\mathfrak{Q}
\circ_{s}^{0\mathbf{Pth}_{\boldsymbol{\mathcal{A}}^{(1)}}}
\mathfrak{P}
\right)
=
\tau^{\mathbf{Pth}_{\boldsymbol{\mathcal{B}}^{(1)}}^{\mathbf{f}^{(1)}}}\left(\left(
f_{s_{j}}^{(1)\flat}\left(
\mathfrak{Q}_{j} \circ_{s_{j}}^{0\mathbf{Pth}_{\boldsymbol{\mathcal{A}}^{(1)}}} \mathfrak{P}_{j}
\right)
\right)_{j\in\bb{\mathbf{s}}}\right).
$$

Finally, note that, for every $j\in\bb{\mathbf{s}}$, $(\mathfrak{Q}_{j} \circ_{s}^{0\mathbf{Pth}_{\boldsymbol{\mathcal{A}}^{(1)}}} \mathfrak{P}_{j}, s_{j}) \prec_{\mathbf{Pth}_{\boldsymbol{\mathcal{A}}^{(1)}}} (\mathfrak{Q} \circ_{s}^{0\mathbf{Pth}_{\boldsymbol{\mathcal{A}}^{(1)}}} \mathfrak{P}, s)$. Thus, by induction, we have that
$$
\left[
f_{s_{j}}^{(1)\flat}\left(
\mathfrak{Q}_{j} \circ_{s_{j}}^{0\mathbf{Pth}_{\boldsymbol{\mathcal{A}}^{(1)}}} \mathfrak{P}_{j}
\right)
\right]_{\varphi(s_{j})}
=
\left[
f_{s_{j}}^{(1)\flat}\left(
\mathfrak{Q}_{j}
\right)
\right]_{\varphi(s_{j})}
\circ_{s_{j}}^{0[\mathbf{Pth}_{\boldsymbol{\mathcal{B}}^{(1)}}^{\mathbf{f}^{(1)}}]}
\left[
f_{s_{j}}^{(1)\flat}\left(
\mathfrak{P}_{j}
\right)
\right]_{\varphi(s_{j})}.
$$

The following chain of equalities holds
\begin{flushleft}
$
\left[
f_{s}^{(1)\flat}\left(
\mathfrak{Q} \circ_{s}^{0\mathbf{Pth}_{\boldsymbol{\mathcal{A}}^{(1)}}} \mathfrak{P}
\right)
\right]_{\varphi(s)}
$
\allowdisplaybreaks
\begin{align*}
&=
\left[
\sigma^{\mathbf{Pth}_{\boldsymbol{\mathcal{B}}^{(1)}}^{\mathbf{f}^{(1)}}}\left(\left(
f_{s_{j}}^{(1)\flat}\left(
\mathfrak{Q}_{j} \circ_{s_{j}}^{0\mathbf{Pth}_{\boldsymbol{\mathcal{A}}^{(1)}}} \mathfrak{P}_{j}
\right)
\right)_{j\in\bb{\mathbf{s}}}\right)
\right]_{\varphi(s)}
\tag{1}
\\
&=
\sigma^{[\mathbf{Pth}_{\boldsymbol{\mathcal{B}}^{(1)}}^{\mathbf{f}^{(1)}}]} \left(\left(
\left[
f_{s_{j}}^{(1)\flat}\left(
\mathfrak{Q}_{j} \circ_{s_{j}}^{0\mathbf{Pth}_{\boldsymbol{\mathcal{A}}^{(1)}}} \mathfrak{P}_{j}
\right)
\right]_{\varphi(s_{j})}
\right)_{j\in\bb{\mathbf{s}}}\right)
\tag{2}
\\
&=
\resizebox{.89\textwidth}{!}{%
$
\sigma^{[\mathbf{Pth}_{\boldsymbol{\mathcal{B}}^{(1)}}^{\mathbf{f}^{(1)}}]} \left(\left(
\left[
f_{s_{j}}^{(1)\flat}\left(
\mathfrak{Q}_{j}
\right)
\right]_{\varphi(s_{j})}
\circ_{s_{j}}^{0[\mathbf{Pth}_{\boldsymbol{\mathcal{B}}^{(1)}}^{\mathbf{f}^{(1)}}]}
\left[
f_{s_{j}}^{(1)\flat}\left(
\mathfrak{P}_{j}
\right)
\right]_{\varphi(s_{j})}
\right)_{j\in\bb{\mathbf{s}}}\right)
$
}
\tag{3}
\\
&=
\sigma^{[\mathbf{Pth}_{\boldsymbol{\mathcal{B}}^{(1)}}^{\mathbf{f}^{(1)}}]} \left(\left(
\left[
f_{s_{j}}^{(1)\flat}\left(
\mathfrak{Q}_{j}
\right)
\right]_{\varphi(s_{j})}
\right)_{j\in\bb{\mathbf{s}}}\right)
\circ_{s}^{0[\mathbf{Pth}_{\boldsymbol{\mathcal{B}}^{(1)}}^{\mathbf{f}^{(1)}}]}
\\
&\hspace{6cm}
\sigma^{[\mathbf{Pth}_{\boldsymbol{\mathcal{B}}^{(1)}}^{\mathbf{f}^{(1)}}]} \left(\left(
\left[
f_{s_{j}}^{(1)\flat}\left(
\mathfrak{P}_{j}
\right)
\right]_{\varphi(s_{j})}
\right)_{j\in\bb{\mathbf{s}}}\right)
\tag{4}
\\
&=
\left[
\sigma^{\mathbf{Pth}_{\boldsymbol{\mathcal{B}}^{(1)}}^{\mathbf{f}^{(1)}}}\left(\left(
f_{s_{j}}^{(1)\flat}\left(
\mathfrak{Q}_{j}
\right)
\right)_{j\in\bb{\mathbf{s}}}\right)
\right]_{\varphi(s)}
\circ_{s}^{0[\mathbf{Pth}_{\boldsymbol{\mathcal{B}}^{(1)}}^{\mathbf{f}^{(1)}}]}
\\
&\hspace{6.3cm}
\left[
\sigma^{\mathbf{Pth}_{\boldsymbol{\mathcal{B}}^{(1)}}^{\mathbf{f}^{(1)}}}\left(\left(
f_{s_{j}}^{(1)\flat}\left(
\mathfrak{P}_{j}
\right)
\right)_{j\in\bb{\mathbf{s}}}\right)
\right]_{\varphi(s)}
\tag{5}
\\
&=
\left[
f_{s}^{(1)\flat}\left(
\mathfrak{Q}
\right)
\right]_{\varphi(s)}
\circ_{s}^{0[\mathbf{Pth}_{\boldsymbol{\mathcal{B}}^{(1)}}^{\mathbf{f}^{(1)}}]}
\left[
f_{s}^{(1)\flat}\left(
\mathfrak{P}
\right)
\right]_{\varphi(s)}.
\tag{6}
\end{align*}
\end{flushleft}

The first equality unravels the definition of the path extension mapping $f^{(1)\flat}$ introduced in Proposition~\ref{PPthExt} at an echelonless path;
the second equality recovers the interpretation of the operation symbol $\sigma$ in the $\Sigma$-algebra $[\mathbf{Pth}_{\boldsymbol{\mathcal{B}}^{(1)}}^{\mathbf{f}^{(1)}(0,1)}]$ introduced in Remark~\ref{RDSigmaAlg};
the third equality follows by Artinian induction;
the fourth equality follows from Proposition~\ref{PQPthBVarA8};
the fifth equality unravels the interpretation of the operation symbol $\sigma$ in the partial $\Sigma$-algebra $[\mathbf{Pth}_{\boldsymbol{\mathcal{B}}^{(1)}}^{\mathbf{f}^{(1)}(0,1)}]$ introduced in Remark~\ref{RDSigmaAlg};
finally, the last equality recovers the definition of $f^{(1)\flat}$ at an echelonless path.

This completes case $(2)$.

Hence, $\mathrm{pr}_{\boldsymbol{\mathcal{B}}^{(1)}, \varphi}^{\mathrm{Ker}(\mathrm{CH}^{(1)})} \circ f^{(1)\flat}$ is compatible with the $1$-composition operation.

All in all, we conclude that $\mathrm{pr}_{\boldsymbol{\mathcal{B}}^{(1)}, \varphi}^{\mathrm{Ker}(\mathrm{CH}^{(1)})} \circ f^{(1)\flat}$ is a $\Sigma^{\boldsymbol{\mathcal{A}}^{(1)}}$-homomorphism.

It remains to prove that
$$
\mathrm{Ker}(\mathrm{CH}^{(1)}_{\boldsymbol{\mathcal{A}}^{(1)}})
\subseteq
\mathrm{Ker}\left( \mathrm{pr}_{\boldsymbol{\mathcal{B}}^{(1)}, \varphi}^{\mathrm{Ker}(\mathrm{CH}^{(1)})} \circ f^{(1)\flat} \right).
$$

Let $s$ be a sort in $S$ and let $\mathfrak{Q}, \mathfrak{P}$ be two paths in $\mathrm{Pth}_{\boldsymbol{\mathcal{A}}^{(1)}, s}$ satisfying that $(\mathfrak{Q}, \mathfrak{P}) \in \mathrm{Ker}(\mathrm{CH}^{(1)}_{\boldsymbol{\mathcal{A}}^{(1)}})_{s}$.

we will prove that
$$
\left[
f^{(1)\flat}_{s}\left(
\mathfrak{P}
\right)
\right]_{\varphi(s)}
=
\left[
f^{(1)\flat}_{s}\left(
\mathfrak{Q}
\right)
\right]_{\varphi(s)}
$$
by Artinian induction on $(\coprod \mathrm{Pth}_{\boldsymbol{\mathcal{A}}^{(1)}}, \leq_{\mathbf{Pth}_{\boldsymbol{\mathcal{A}}^{(1)}}})$.

{\sffamily Base step of the Artinian induction}

Let $(\mathfrak{P}, s)$ be a minimal element of $(\coprod \mathrm{Pth}_{\boldsymbol{\mathcal{A}}^{(1)}}, \leq_{\mathbf{Pth}_{\boldsymbol{\mathcal{A}}^{(1)}}})$. Then, by Proposition~\ref{PMinimal}, the path $\mathfrak{P}$ is either (1) a $(1,0)$-identity path or (2) an echelon.

In either case, the equality
$$
\left[
f^{(1)\flat}_{s}\left(
\mathfrak{P}
\right)
\right]_{\varphi(s)}
=
\left[
f^{(1)\flat}_{s}\left(
\mathfrak{Q}
\right)
\right]_{\varphi(s)}
$$
follows  directly from the fact that, according to Corollary~\ref{CCHUZId} or Proposition~\ref{PCHEch}, $\mathfrak{P} = \mathfrak{Q}$.

{\sffamily Inductive step of the Artinian induction}

Let $(\mathfrak{P}, s)$ be a non-minimal element of $(\coprod\mathrm{Pth}_{\boldsymbol{\mathcal{A}}^{(1)}}, \leq_{\mathbf{Pth}_{\boldsymbol{\mathcal{A}}^{(1)}}})$. Let us suppose that, for every sort $t \in S$ and every path $\mathfrak{P}'$ in $\mathrm{Pth}_{\boldsymbol{\mathcal{A}}^{(1)}, t}$, if $(\mathfrak{P}',t) <_{\mathbf{Pth}_{\boldsymbol{\mathcal{A}}^{(1)}}} (\mathfrak{P}, s)$, the the statement holds for $\mathfrak{P}'$, i.e., for every path $\mathfrak{Q}'$ in $\mathrm{Pth}_{\boldsymbol{\mathcal{A}}^{(1)}, t}$, if $(\mathfrak{P}', \mathfrak{Q}')\in\mathrm{Ker}(\mathrm{CH}^{(1)})_{t}$, then 
$$
\left[
f^{(1)\flat}_{t}\left(
\mathfrak{P}'
\right)
\right]_{\varphi(t)}
=
\left[
f^{(1)\flat}_{t}\left(
\mathfrak{Q}'
\right)
\right]_{\varphi(t)}.
$$

Let $(\mathfrak{P}, s)$ be a non-minimal element of $(\coprod \mathrm{Pth}_{\boldsymbol{\mathcal{A}}^{(1)}, \leq_{\mathbf{Pth}_{\boldsymbol{\mathcal{A}}^{(1)}}}})$. We can assume that $\mathfrak{P}$ is not a $(1,0)$-identity path, since that case has already been considered. By Lemma~\ref{LDOrdI}, that $\mathfrak{P}$ is either (1) a path of length strictly greater than one containing at least one echelon or (2) an echelonless path.

If (1), then let $i \in \bb{\mathfrak{P}}$ be the first index for which the one-step subpath $\mathfrak{P}^{i,i}$ is an echelon. We distinguish two cases accordingly.

If $i=0$, i.e., if $\mathfrak{P}$ has its first echelon on its first step, then according to Proposition~\ref{PPthExt}, we have that
$$
f^{(1)\flat}_{s}\left(
\mathfrak{P}
\right)
=
f^{(1)\flat}_{s}\left(
\mathfrak{P}^{1,\bb{\mathfrak{P}}-1}
\right)
\circ_{\varphi(s)}^{0\mathbf{Pth}_{\boldsymbol{\mathcal{B}}^{(1)}}}
f^{(1)\flat}_{s}\left(
\mathfrak{P}^{0,0}
\right).
$$

Since $\mathrm{CH}^{(1)}_{\boldsymbol{\mathcal{A}}^{(1)}, s}(\mathfrak{P}) \in \eta^{(1,\mathcal{A}^{(1)})}[\mathcal{A}^{(1)}]_{s}^{\mathrm{int}}$ and $(\mathfrak{P}, \mathfrak{Q}) \in \mathrm{Ker}(\mathrm{CH}^{(1)}_{\boldsymbol{\mathcal{A}}^{(1)}})_{s}$, we have, by Lemma~\ref{LCHEchInt}, that $\mathfrak{Q}$ is a path of length strictly greater than one containing its first echelon on its first step.

Thus, according to Proposition~\ref{PPthExt}, we have that
$$
f^{(1)\flat}_{s}\left(
\mathfrak{Q}
\right)
=
f^{(1)\flat}_{s}\left(
\mathfrak{Q}^{1,\bb{\mathfrak{P}}-1}
\right)
\circ_{\varphi(s)}^{0\mathbf{Pth}_{\boldsymbol{\mathcal{B}}^{(1)}}}
f^{(1)\flat}_{s}\left(
\mathfrak{Q}^{0,0}
\right).
$$

Since $(\mathfrak{P}, \mathfrak{Q}) \in \mathrm{Ker}(\mathrm{CH}^{(1)}_{\boldsymbol{\mathcal{A}}^{(1)}})_{s}$, we have that the pairs $(\mathfrak{P}^{1,\bb{\mathfrak{P}}-1}, \mathfrak{Q}^{1,\bb{\mathfrak{Q}}-1}$ and $(\mathfrak{P}^{0,0}, \mathfrak{Q}^{0,0})$ are in $\mathrm{Ker}(\mathrm{CH}^{(1)})_{\boldsymbol{\mathcal{A}}^{(1)}, s}$. Note that, according to Definition~\ref{DOrd}, we have that $(\mathfrak{P}^{0,0}, s)$ and $(\mathfrak{P}^{1,\bb{\mathfrak{P}}-1}, s)$$\prec_{\mathbf{Pth}_{\boldsymbol{\mathcal{A}}^{(1)}}}$-precede $(\mathfrak{P}, s)$.

Therefore, by the inductive hypothesis, the paths $\mathfrak{P}^{0,0}$ and $\mathfrak{Q}^{0,0}$, and the paths $\mathfrak{P}^{1,\bb{\mathfrak{P}}-1}$ and $\mathfrak{Q}^{1,\bb{\mathfrak{Q}}-1}$ satisfy
\begin{align*}
\left[
f^{(1)\flat}_{s}\left(
\mathfrak{P}^{0,0}
\right)
\right]_{\varphi(s)}
&=
\left[
f^{(1)\flat}_{s}\left(
\mathfrak{Q}^{0,0}
\right)
\right]_{\varphi(s)}
\\
\left[
f^{(1)\flat}_{s}\left(
\mathfrak{P}^{1,\bb{\mathfrak{P}}-1}
\right)
\right]_{\varphi(s)}
&=
\left[
f^{(1)\flat}_{s}\left(
\mathfrak{Q}^{1,\bb{\mathfrak{Q}}-1}
\right)
\right]_{\varphi(s)}.
\end{align*}

Thus, the following chain of equalities holds
\allowdisplaybreaks
\begin{align*}
\left[
f^{(1)\flat}_{s}\left(
\mathfrak{P}
\right)
\right]_{\varphi(s)}
&=
\left[
f^{(1)\flat}_{s}\left(
\mathfrak{P}^{1,\bb{\mathfrak{P}}-1}
\right)
\circ_{\varphi(s)}^{0\mathbf{Pth}_{\boldsymbol{\mathcal{B}}^{(1)}}}
f^{(1)\flat}_{s}\left(
\mathfrak{P}^{0,0}
\right)
\right]_{\varphi(s)}
\tag{1}
\\
&=
\left[
f^{(1)\flat}_{s}\left(
\mathfrak{P}^{1,\bb{\mathfrak{P}}-1}
\right)
\right]_{\varphi(s)}
\circ_{\varphi(s)}^{0[\mathbf{Pth}_{\boldsymbol{\mathcal{B}}^{(1)}}]}
\left[
f^{(1)\flat}_{s}\left(
\mathfrak{P}^{0,0}
\right)
\right]_{\varphi(s)}
\tag{2}
\\
&=
\left[
f^{(1)\flat}_{s}\left(
\mathfrak{Q}^{1,\bb{\mathfrak{Q}}-1}
\right)
\right]_{\varphi(s)}
\circ_{\varphi(s)}^{0[\mathbf{Pth}_{\boldsymbol{\mathcal{B}}^{(1)}}]}
\left[
f^{(1)\flat}_{s}\left(
\mathfrak{Q}^{0,0}
\right)
\right]_{\varphi(s)}
\tag{3}
\\
&=
\left[
f^{(1)\flat}_{s}\left(
\mathfrak{Q}^{1,\bb{\mathfrak{Q}}-1}
\right)
\circ_{\varphi(s)}^{0\mathbf{Pth}_{\boldsymbol{\mathcal{B}}^{(1)}}}
f^{(1)\flat}_{s}\left(
\mathfrak{Q}^{0,0}
\right)
\right]_{\varphi(s)}
\tag{4}
\\
&=
\left[
f^{(1)\flat}_{s}\left(
\mathfrak{Q}
\right)
\right]_{\varphi(s)}.
\tag{5}
\end{align*}

The case of $\mathfrak{P}$ being a path of length strictly greater than one containing its first echelon on its first step follows.

If $i>0$, that is, if $\mathfrak{P}$ is a path of length strictly greater than one containing its first echelon on a step different from the initial one, then according to Proposition~\ref{PPthExt}, we have that
$$
f^{(1)\flat}_{s}\left(
\mathfrak{P}
\right)
=
f^{(1)\flat}_{s}\left(
\mathfrak{P}^{i,\bb{\mathfrak{P}}-1}
\right)
\circ_{\varphi(s)}^{0\mathbf{Pth}_{\boldsymbol{\mathcal{B}}^{(1)}}}
f^{(1)\flat}_{s}\left(
\mathfrak{P}^{0,i-1}
\right).
$$

Since $\mathrm{CH}^{(1)}_{\boldsymbol{\mathcal{A}}^{(1)}, s}(\mathfrak{P}) \in \eta^{(1,\mathcal{A}^{(1)})}[\mathcal{A}^{(1)}]_{s}^{\neg\mathrm{int}}$ and $(\mathfrak{P}, \mathfrak{Q}) \in \mathrm{Ker}(\mathrm{CH}^{(1)}_{\boldsymbol{\mathcal{A}}^{(1)}})_{s}$, we have, by Lemma~\ref{LCHEchNInt}, that $\mathfrak{Q}$ is a path of length strictly greater than one containing its first echelon on a step different from the initial one.

Thus, according to Proposition~\ref{PPthExt}, we have that
$$
f^{(1)\flat}_{s}\left(
\mathfrak{Q}
\right)
=
f^{(1)\flat}_{s}\left(
\mathfrak{Q}^{i,\bb{\mathfrak{P}}-1}
\right)
\circ_{\varphi(s)}^{0\mathbf{Pth}_{\boldsymbol{\mathcal{B}}^{(1)}}}
f^{(1)\flat}_{s}\left(
\mathfrak{Q}^{0,i-1}
\right).
$$

Since $(\mathfrak{P}, \mathfrak{Q}) \in \mathrm{Ker}(\mathrm{CH}^{(1)}_{\boldsymbol{\mathcal{A}}^{(1)}})_{s}$, we have that the pairs $(\mathfrak{P}^{i,\bb{\mathfrak{P}}-1}, \mathfrak{Q}^{i,\bb{\mathfrak{Q}}-1}$ and $(\mathfrak{P}^{0,i-1}, \mathfrak{Q}^{0,i-1})$ are in $\mathrm{Ker}(\mathrm{CH}^{(1)})_{s}$. Note that, according to Definition~\ref{DOrd}, we have that $(\mathfrak{P}^{0,i-1}, s)$ and $(\mathfrak{P}^{i,\bb{\mathfrak{P}}-1}, s)$$\prec_{\mathbf{Pth}_{\boldsymbol{\mathcal{A}}^{(1)}}}$-precede $(\mathfrak{P}, s)$.

Therefore, by the inductive hypothesis, the paths $\mathfrak{P}^{0,i-1}$ and $\mathfrak{Q}^{0,i-1}$, and the paths $\mathfrak{P}^{i,\bb{\mathfrak{P}}-1}$ and $\mathfrak{Q}^{i,\bb{\mathfrak{Q}}-1}$ satisfy
\begin{align*}
\left[
f^{(1)\flat}_{s}\left(
\mathfrak{P}^{0,i-1}
\right)
\right]_{\varphi(s)}
&=
\left[
f^{(1)\flat}_{s}\left(
\mathfrak{Q}^{0,i-1}
\right)
\right]_{\varphi(s)}
\\
\left[
f^{(1)\flat}_{s}\left(
\mathfrak{P}^{i,\bb{\mathfrak{P}}-1}
\right)
\right]_{\varphi(s)}
&=
\left[
f^{(1)\flat}_{s}\left(
\mathfrak{Q}^{i,\bb{\mathfrak{Q}}-1}
\right)
\right]_{\varphi(s)}.
\end{align*}

Thus, the following chain of equalities holds
\allowdisplaybreaks
\begin{align*}
\left[
f^{(1)\flat}_{s}\left(
\mathfrak{P}
\right)
\right]_{\varphi(s)}
&=
\left[
f^{(1)\flat}_{s}\left(
\mathfrak{P}^{i,\bb{\mathfrak{P}}-1}
\right)
\circ_{\varphi(s)}^{0\mathbf{Pth}_{\boldsymbol{\mathcal{B}}^{(1)}}}
f^{(1)\flat}_{s}\left(
\mathfrak{P}^{0,i-1}
\right)
\right]_{\varphi(s)}
\tag{1}
\\
&=
\left[
f^{(1)\flat}_{s}\left(
\mathfrak{P}^{i,\bb{\mathfrak{P}}-1}
\right)
\right]_{\varphi(s)}
\circ_{\varphi(s)}^{0[\mathbf{Pth}_{\boldsymbol{\mathcal{B}}^{(1)}}]}
\left[
f^{(1)\flat}_{s}\left(
\mathfrak{P}^{0,i-1}
\right)
\right]_{\varphi(s)}
\tag{2}
\\
&=
\left[
f^{(1)\flat}_{s}\left(
\mathfrak{Q}^{i,\bb{\mathfrak{Q}}-1}
\right)
\right]_{\varphi(s)}
\circ_{\varphi(s)}^{0[\mathbf{Pth}_{\boldsymbol{\mathcal{B}}^{(1)}}]}
\left[
f^{(1)\flat}_{s}\left(
\mathfrak{Q}^{0,i-1}
\right)
\right]_{\varphi(s)}
\tag{3}
\\
&=
\left[
f^{(1)\flat}_{s}\left(
\mathfrak{Q}^{i,\bb{\mathfrak{Q}}-1}
\right)
\circ_{\varphi(s)}^{0\mathbf{Pth}_{\boldsymbol{\mathcal{B}}^{(1)}}}
f^{(1)\flat}_{s}\left(
\mathfrak{Q}^{0,i-1}
\right)
\right]_{\varphi(s)}
\tag{4}
\\
&=
\left[
f^{(1)\flat}_{s}\left(
\mathfrak{Q}
\right)
\right]_{\varphi(s)}
\tag{5}
\end{align*}

The case of $\mathfrak{P}$ being a path of length strictly greater than one containing its first echelon on a step different from the initial one follows.

Case (1) follows.

If (2), i.e., if $\mathfrak{P}$ is an echelonless path, then there exists a unique word $\mathbf{s} \in S^{\star} - \{\lambda\}$ and a unique operation symbol $\sigma \in \Sigma_{\mathbf{s}, s}$ associated to $\mathfrak{P}$. Let $(\mathfrak{P}_{j})_{j \in \bb{\mathbf{s}}}$ be the family of paths in $\mathbf{Pth}_{\boldsymbol{\mathcal{A}}^{(1)}, \mathbf{s}}$ which, in virtue of Lemma~\ref{LPthExtract}, we can extract from $\mathfrak{P}$. Then, according to Proposition~\ref{PPthExt}, we have that
$$
f^{(1)\flat}_{s}\left(
\mathfrak{P}
\right)
=
\tau^{\mathbf{Pth}_{\boldsymbol{\mathcal{B}}^{(1)}}^{\mathbf{f}^{(1)}}}\left(\left(
f^{(1)\flat}_{s_{j}}\left(
\mathfrak{P}_{j}
\right)
\right)_{j\in\bb{\mathbf{s}}}\right).
$$

Since $\mathrm{CH}^{(1)}_{\boldsymbol{\mathcal{A}}^{(1)}, s}(\mathfrak{P}) \in \mathcal{T}(\sigma, \mathrm{T}_{\Sigma^{\boldsymbol{\mathcal{A}}^{(1)}}(X)})_{1}$, which is a subset of $\mathrm{T}_{\Sigma^{\boldsymbol{\mathcal{A}}^{(1)}}}(X)_{s}^{\mathsf{E}}$, and $(\mathfrak{P}, \mathfrak{Q}) \in \mathrm{Ker}(\mathrm{CH}^{(1)}_{\boldsymbol{\mathcal{A}}^{(1)}})_{s}$ we have, by Lemma~\ref{LCHNEch}, that $\mathfrak{Q}$ is an echelonless path associated to $\tau$, the same operation symbol as that associated to $\mathfrak{P}$.

Let $(\mathfrak{Q}_{j})_{j\in\bb{\mathbf{s}}}$ be the family of paths in $\mathrm{Pth}_{\boldsymbol{\mathcal{A}}^{(1)}, \mathbf{s}}$ which, by Lemma~\ref{LPthExtract}, we can extract from $\mathfrak{Q}$. Then, according to Proposition~\ref{PPthExt}, we have that 
$$
f^{(1)\flat}_{s}\left(
\mathfrak{Q}
\right)
=
\sigma^{\mathbf{Pth}_{\boldsymbol{\mathcal{B}}^{(1)}}^{\mathbf{f}^{(1)}(0,1)}}\left(\left(
f^{(1)\flat}_{s_{j}}\left(
\mathfrak{Q}_{j}
\right)
\right)_{j\in\bb{\mathbf{s}}}\right).
$$

Since $(\mathfrak{P}, \mathfrak{Q}) \in \mathrm{Ker}(\mathrm{CH}^{(1)}_{\boldsymbol{\mathcal{A}}^{(1)}})_{s}$, we have, for every $j \in \bb{\mathbf{s}}$, that $(\mathfrak{P}_{j}, \mathfrak{Q}_{j}) \in \mathrm{Ker}(\mathrm{CH}^{(1)}_{\boldsymbol{\mathcal{A}}^{(1)}})_{s_{j}}$. Note that, according to Definition~\ref{DOrd}, we have that, for every $j \in \bb{\mathbf{s}}$, $(\mathfrak{P}_{j}, s_{j})$$\prec_{\mathbf{Pth}_{\boldsymbol{\mathcal{A}}^{(1)}}}$-precedes $(\mathfrak{P},s)$.

Therefore, by the inductive hypothesis, for every $j \in \bb{\mathbf{s}}$, the paths $\mathfrak{P}_{j}$ and $\mathfrak{Q}_{j}$ satisfy
$$
\left[
f^{(1)\flat}_{s_{j}}\left(
\mathfrak{P}_{j}
\right)
\right]_{\varphi(s_{j})}
=
\left[
f^{(1)\flat}_{s_{j}}\left(
\mathfrak{Q}_{j}
\right)
\right]_{\varphi(s_{j})}
$$

Thus, the following chain of equalities holds
\allowdisplaybreaks
\begin{align*}
\left[
f^{(1)\flat}_{s}\left(
\mathfrak{P}
\right)
\right]_{\varphi(s)}
&=
\left[
\sigma^{\mathbf{Pth}_{\boldsymbol{\mathcal{B}}^{(1)}}^{\mathbf{f}^{(1)}}}\left(\left(
f^{(1)\flat}_{s_{j}}\left(
\mathfrak{P}_{j}
\right)
\right)_{j\in\bb{\mathbf{s}}}\right)
\right]_{\varphi(s)}
\tag{1}
\\
&=
\sigma^{[\mathbf{Pth}_{\boldsymbol{\mathcal{B}}^{(1)}}^{\mathbf{f}^{(1)}(0,1)}]}\left(\left(
\left[
f^{(1)\flat}_{s_{j}}\left(
\mathfrak{P}_{j}
\right)
\right]_{\varphi(s_{j})}
\right)_{j\in\bb{\mathbf{s}}}\right)
\tag{2}
\\
&=
\sigma^{[\mathbf{Pth}_{\boldsymbol{\mathcal{B}}^{(1)}}^{\mathbf{f}^{(1)}(0,1)}]}\left(\left(
\left[
f^{(1)\flat}_{s_{j}}\left(
\mathfrak{Q}_{j}
\right)
\right]_{\varphi(s_{j})}
\right)_{j\in\bb{\mathbf{s}}}\right)
\tag{3}
\\
&=
\left[
\sigma^{\mathbf{Pth}_{\boldsymbol{\mathcal{B}}^{(1)}}^{\mathbf{f}^{(1)}(0,1)}}\left(\left(
f^{(1)\flat}_{s_{j}}\left(
\mathfrak{Q}_{j}
\right)
\right)_{j\in\bb{\mathbf{s}}}\right)
\right]_{\varphi(s)}
\tag{4}
\\
&=
\left[
f^{(1)\flat}_{s}\left(
\mathfrak{Q}
\right)
\right]_{\varphi(s)}
\tag{5}
\end{align*}

The case of $\mathfrak{P}$ being an echelonless path follows.

Case (1) follows.

This proves that $\mathrm{Ker}(\mathrm{CH}^{(1)}_{\boldsymbol{\mathcal{A}}^{(1)}})$ is included in $\mathrm{Ker}(\mathrm{pr}_{\boldsymbol{\mathcal{B}}^{(1)}, \varphi}^{\mathrm{Ker}(\mathrm{CH}^{(1)})} \circ f^{(1)\flat})$.

This completes the proof of Proposition~\ref{PHomPthExtKer}.
\end{proof}

Thus, according to the universal property of the quotient, we introduce the following definition.

\begin{definition}
\label{DQPthExt}
Following Proposition~\ref{PHomPthExtKer} and taking into account the Universal Property of the Quotient, there exists a unique $\Sigma^{\boldsymbol{\mathcal{A}}^{(1)}}$-homomorphism, that we will denote by $f^{[1] @}$, i.e., 
$$
f^{[1] @}
\colon
[ \mathbf{Pth}_{\boldsymbol{\mathcal{A}^{(1)}}}]
\mor
[ \mathbf{Pth}_{\boldsymbol{\mathcal{B}}^{(1)}}^{\mathbf{f}^{(1)}}]
$$
satisfying that $f^{[1] @} \circ \mathrm{pr}^{\mathrm{Ker}(\mathrm{CH}^{(1)})}_{\boldsymbol{\mathcal{A}}^{(1)}} = \mathrm{pr}_{\boldsymbol{\mathcal{B}}^{(1)}, \varphi}^{\mathrm{Ker}(\mathrm{CH}^{(1)})} \circ f^{(1)\flat}$, namely 
$$
f^{[1] @}
=
\left(
\mathrm{pr}_{\boldsymbol{\mathcal{B}}^{(1)}, \varphi}^{\mathrm{Ker}(\mathrm{CH}^{(1)})} \circ f^{(1)\flat}
\right)^{\natural}.
$$
We will call this mapping the \emph{first-order quotient path extension mapping} of $\mathbf{f}^{(1)}$. Formally, for every sort $s$ in $S$ and every path class $[ \mathfrak{P} ]_{s}$ in $[ \mathrm{Pth}_{\boldsymbol{\mathcal{A}}^{(1)}} ]_{s}$, 
$$
f^{[1] @}_{s}\left(
\left[
\mathfrak{P}
\right]_{s}
\right)
=
\left[
f^{(1)\flat}_{s}\left(
\mathfrak{P}
\right)
\right]_{\varphi(s)}
$$

Recall that, taking into account Theorems~\ref{TIso} and \ref{TPthBPTB},  we obtain a unique $\Sigma^{\boldsymbol{\mathcal{A}}^{(1)}}$-homomorphism from $[\mathbf{PT}_{\boldsymbol{\mathcal{A}}^{(1)}}]$ to $[\mathbf{PT}_{\boldsymbol{\mathcal{B}}^{(1)}}^{\mathbf{f}^{(1)}}]$. We agreed to also use $f^{[1] @}$ to denote this mapping
$$
f^{[1] @}
\colon
[ \mathbf{PT}_{\boldsymbol{\mathcal{A}^{(1)}}}]
\mor
[ \mathbf{PT}_{\boldsymbol{\mathcal{B}}^{(1)}}^{\mathbf{f}^{(1)}}].
$$
Formally, for every sort $s$ in $S$ and every path term class $[ P ]_{s}$ in $[ \mathrm{PT}_{\boldsymbol{\mathcal{A}}^{(1)}}]_{s}$,
$$
f^{[1] @}_{s}\left(
\left[
P
\right]_{s}
\right)
=
\left[
\mathrm{CH}^{(1)}_{\boldsymbol{\mathcal{B}}^{(1)}, \varphi(s)}\left(
f^{(1)\flat}_{s}\left(
\mathrm{ip}^{(1,X)@}_{\boldsymbol{\mathcal{A}}^{(1)}, s}\left(
P
\right)
\right)
\right)
\right]_{\varphi(s)}.
$$

\end{definition}

\section{The behaviour of the first-order quotient path extension mapping}

To end this chapter, we study the relations between the first-order quotient path-extension mapping and the mappings $\mathrm{CH}^{[ 1 ]}$ and $\mathrm{ip}^{([ 1 ], X)@}$ and the source, target and identity path mappings.

\begin{proposition}
\label{PQPthExtQCH}
Let $\mathbf{f}^{(1)}=(\varphi, c, (f^{(i)})_{i\in 2})$ be a morphism of first-order many-sorted rewriting systems from $\boldsymbol{\mathcal{A}}^{(1)}$ to $\boldsymbol{\mathcal{B}}^{(1)}$. Then the following equalities holds
$$
f^{[1] @} \circ \mathrm{CH}_{\boldsymbol{\mathcal{A}}^{(1)}}^{[1]}
=
\mathrm{CH}_{\boldsymbol{\mathcal{B}}^{(1)}, \varphi}^{[1]} \circ f^{[1] @}.
$$
\end{proposition}

\begin{proof}
For every sort $s$ in $S$ and every path class $[ \mathfrak{P} ]_{s}$ in $[ \mathrm{Pth}_{\boldsymbol{\mathcal{A}}^{(1)}} ]_{s}$, the following chain of equalities holds
\begin{flushleft}
$
f^{[ 1 ] @}_{s} \left(
\mathrm{CH}^{[ 1 ]}_{\boldsymbol{\mathcal{A}}^{(1)}, s} \left(
\left[
\mathfrak{P}
\right]_{s}
\right)
\right)
$
\allowdisplaybreaks
\begin{align*}
&=
f^{[ 1 ] @}_{s} \left(
\left[
\mathrm{CH}^{(1)}_{\boldsymbol{\mathcal{A}}^{(1)}, s} \left(
\mathfrak{P}
\right)
\right]_{s}
\right)
\tag{1}
\\
&=
\left[
\mathrm{CH}^{(1)}_{\boldsymbol{\mathcal{B}}^{(1)}, \varphi(s)} \left(
f^{(1)\flat}_{s} \left(
\mathrm{ip}^{(1,X)@}_{\boldsymbol{\mathcal{A}}^{(1)}, s}\left(
\mathrm{CH}^{(1)}_{\boldsymbol{\mathcal{A}}^{(1)}, s} \left(
\mathfrak{P}
\right)
\right)
\right)
\right)
\right]_{\varphi(s)}
\tag{2}
\\
&=
\mathrm{CH}^{[ 1 ]}_{\boldsymbol{\mathcal{B}}^{(1)}, \varphi(s)} \left(
\left[
f^{(1) \flat}_{s} \left(
\mathrm{ip}^{(1,X)@}_{\boldsymbol{\mathcal{A}}^{(1)}, s}\left(
\mathrm{CH}^{(1)}_{\boldsymbol{\mathcal{A}}^{(1)}, s} \left(
\mathfrak{P}
\right)
\right)
\right)
\right]_{\varphi(s)}
\right)
\tag{3}
\\
&=
\mathrm{CH}^{[ 1 ]}_{\boldsymbol{\mathcal{B}}^{(1)}, \varphi(s)} \left(
f^{[1] @}_{s} \left(
\left[
\mathrm{ip}^{(1,X)@}_{\boldsymbol{\mathcal{A}}^{(1)}, s}\left(
\mathrm{CH}^{(1)}_{\boldsymbol{\mathcal{A}}^{(1)}, s} \left(
\mathfrak{P}
\right)
\right)
\right]_{s}
\right)
\right)
\tag{4}
\\
&=
\mathrm{CH}^{[ 1 ]}_{\boldsymbol{\mathcal{B}}^{(1)}, \varphi(s)} \left(
f^{[1] @}_{s} \left(
\left[
\mathfrak{P}
\right]_{s}
\right)
\right)
\tag{5}
\end{align*}
\end{flushleft}

The first equality unravels the definition of $s$-th component of the $S$-sorted mapping $\mathrm{CH}^{[ 1 ]}_{\boldsymbol{\mathcal{A}}^{(1)}}$, introduced in Definition~\ref{DPTQCH};
the second equality unravels the definition of the $s$-th component of the $S$-sorted mapping $f^{[1] @}$ at a path term class, introduced in Definition~\ref{DQPthExt};
the third equality recovers the definition of the $\varphi(s)$-th component of the $S$-sorted mapping $\mathrm{CH}^{[ 1 ]}_{\boldsymbol{\mathcal{B}}^{(1)}}$, introduced in Definition~\ref{DPTQCH};
the fourth equality recovers the definition of the $s$-th component of the $S$-sorted mapping $f^{[1] @}$ at a path term class, introduced in Definition~\ref{DQPthExt};
finally, the last equality follows from Proposition~\ref{PIpCH}.
\end{proof}

\begin{proposition}
\label{PQPthExtQIp}
Let $\mathbf{f}^{(1)}=(\varphi, c, (f^{(i)})_{i\in 2})$ be a morphism of first-order many-sorted rewriting systems from $\boldsymbol{\mathcal{A}}^{(1)}$ to $\boldsymbol{\mathcal{B}}^{(1)}$. Then the following equalities holds
$$
\mathrm{ip}_{\boldsymbol{\mathcal{B}}^{(1)}, \varphi}^{([1], Y)@}
\circ
f^{[1] @}
=
f^{[1] @}
\circ
\mathrm{ip}_{\boldsymbol{\mathcal{A}}^{(1)}}^{([1], X)@}.
$$
\end{proposition}

\begin{proof}
For every sort $s$ in $S$ and every path term class $[ P ]_{s}$ in $[ \mathrm{PT}_{\boldsymbol{\mathcal{A}}^{(1)}} ]_{s}$, the following chain of equalities holds
\begin{flushleft}
$
\mathrm{ip}_{\boldsymbol{\mathcal{B}}^{(1)}, \varphi(s)}^{([1], Y)@} \left(
f_{s}^{[1] @} \left(
\left[
P
\right]_{s}
\right)
\right)
$
\allowdisplaybreaks
\begin{align*}
&=
\mathrm{ip}_{\boldsymbol{\mathcal{B}}^{(1)}, \varphi(s)}^{([1], Y)@} \left(
\left[
\mathrm{CH}^{(1)}_{\boldsymbol{\mathcal{B}}^{(1)}, \varphi(s)}\left(
f_{s}^{(1) \flat} \left(
\mathrm{ip}^{(1,X)@}_{\boldsymbol{\mathcal{A}}^{(1)}, s}\left(
P
\right)
\right)
\right)
\right]_{\varphi(s)}
\right)
\tag{1}
\\
&=
\left[
\mathrm{ip}_{\boldsymbol{\mathcal{B}}^{(1)}, \varphi(s)}^{(1, Y)@} \left(
\mathrm{CH}^{(1)}_{\boldsymbol{\mathcal{B}}^{(1)}, \varphi(s)}\left(
f_{s}^{(1) \flat} \left(
\mathrm{ip}^{(1,X)@}_{\boldsymbol{\mathcal{A}}^{(1)}, s}\left(
P
\right)
\right)
\right)
\right)
\right]_{\varphi(s)}
\tag{2}
\\
&=
\left[
f_{s}^{(1) \flat} \left(
\mathrm{ip}^{(1,X)@}_{\boldsymbol{\mathcal{A}}^{(1)}, s}\left(
P
\right)
\right)
\right]_{\varphi(s)}
\tag{3}
\\
&=
f_{s}^{[ 1 ] @} \left(
\left[
\mathrm{ip}^{(1,X)@}_{\boldsymbol{\mathcal{A}}^{(1)}, s}\left(
P
\right)
\right]_{s}
\right)
\tag{4}
\\
&=
f_{s}^{[ 1 ] @} \left(
\mathrm{ip}^{(1,X)@}_{\boldsymbol{\mathcal{A}}^{(1)}, s}\left(
\left[
P
\right]_{s}
\right)
\right)
\tag{5}
\end{align*}
\end{flushleft}

The first equality unravels the definition of the $s$-th component of the $S$-sorted mapping $f^{[1] @}$ at a path term class, introduced in Definition~\ref{DQPthExt};
the second equality unravels the definition of the $\varphi(s)$-th component of the $S$-sorted mapping $\mathrm{ip}^{([ 1 ],Y)@}_{\boldsymbol{\mathcal{B}}^{(1)}}$, introduced in Definition~\ref{DPTQIp};
the third equality follows from Proposition~\ref{PIpCH};
the fourth equality recovers the definition of the $s$-th component of the $S$-sorted mapping $f^{[1] @}$ at a path term class, introduced in Definition~\ref{DQPthExt};
finally, the last equality recovers the definition of the $s$-th component of the $S$-sorted mapping $\mathrm{ip}^{([ 1 ],X)@}_{\boldsymbol{\mathcal{A}}^{(1)}}$, introduced in Definition~\ref{DPTQIp}.
\end{proof}

\begin{proposition}
\label{PQPthExtScTg}
Let $\mathbf{f}^{(1)}=(\varphi, c, (f^{(i)})_{i\in 2})$ be a morphism of first-order many-sorted rewriting systems from $\boldsymbol{\mathcal{A}}^{(1)}$ to $\boldsymbol{\mathcal{B}}^{(1)}$. Then the following equalities holds
\begin{align*}
\mathrm{sc}_{\boldsymbol{\mathcal{B}}^{(1)},\varphi}^{(0,[1])} \circ f^{[1] @}
&=
f^{(0)\sharp}\circ\mathrm{sc}_{\boldsymbol{\mathcal{A}}^{(1)}}^{(0,[1])}
&&\mbox{and}&
\mathrm{tg}_{\boldsymbol{\mathcal{B}}^{(1)},\varphi}^{(0,[1])} \circ f^{[1] @}
&=
f^{(0)\sharp}\circ\mathrm{tg}_{\boldsymbol{\mathcal{A}}^{(1)}}^{(0,[1])}.
\end{align*}
\end{proposition}

\begin{proof}
For every sort $s$ in $S$ and every path class $[ \mathfrak{P} ]_{s}$ in $[ \mathrm{Pth}_{\boldsymbol{\mathcal{A}}^{(1)}} ]_{s}$, the following chain of equalities holds
\allowdisplaybreaks
\begin{align*}
\mathrm{sc}^{(0,[1])}_{\boldsymbol{\mathcal{B}}^{(1)}, \varphi(s)} \left(
f^{[ 1 ] @}_{s} \left(
\left[ \mathfrak{P} \right]_{s}
\right)
\right)
&=
\mathrm{sc}^{(0,[1])}_{\boldsymbol{\mathcal{B}}^{(1)}, \varphi(s)} \left(
\left[
f^{(1)\flat}_{s} \left(
\mathfrak{P}
\right)
\right]_{\varphi(s)}
\right)
\tag{1}
\\
&=
\mathrm{sc}^{(0,1)}_{\boldsymbol{\mathcal{B}}^{(1)}, \varphi(s)} \left(
f^{(1)\flat}_{s} \left(
\mathfrak{P}
\right)
\right)
\tag{2}
\\
&=
f^{(0)\sharp}_{s} \left(
\mathrm{sc}^{(0,1)}_{\boldsymbol{\mathcal{A}}^{(1)}, s} \left(
\mathfrak{P}
\right)
\right)
\tag{3}
\\
&=
f^{(0)\sharp}_{s} \left(
\mathrm{sc}^{(0,[1])}_{\boldsymbol{\mathcal{A}}^{(1)}, s} \left(
\left[
\mathfrak{P}
\right]_{s}
\right)
\right)
\tag{4}
\end{align*}

The first equality unravels the definition of the first-order quotient path extension mapping $f^{[ 1 ] @}$ at a path class, introduced in Definition~\ref{DQPthExt};
the second equality unravels the definition of the mapping $\mathrm{sc}^{(0,[1])}_{\boldsymbol{\mathcal{B}}^{(1)}, \varphi}$, introduced in Definition~\ref{DCHUZ};
the third equality follows from Proposition~\ref{PPthExt};
finally, the last equality recovers the definition of the mapping $\mathrm{sc}^{(0,[1])}_{\boldsymbol{\mathcal{A}}^{(1)}}$, introduced in Definition~\ref{DCHUZ}.

A similar argument applies to the $(0,[1])$-target mapping.
\end{proof}

\begin{proposition}
\label{PQPthExtIp}
Let $\mathbf{f}^{(1)}=(\varphi, c, (f^{(i)})_{i\in 2})$ be a morphism of first-order many-sorted rewriting systems from $\boldsymbol{\mathcal{A}}^{(1)}$ to $\boldsymbol{\mathcal{B}}^{(1)}$. Then the following equality holds
$$
f^{[1] @} \circ \mathrm{ip}_{\boldsymbol{\mathcal{A}}^{(1)}}^{([1], 0)\sharp}
=
\mathrm{ip}_{\boldsymbol{\mathcal{B}}^{(1)}, \varphi}^{([1],0)\sharp} \circ f^{(0)\sharp}.
$$
\end{proposition}

\begin{proof}
For every sort $s$ in $S$ and every term $P$ in $\mathrm{T}_{\Sigma}(X)_{s}$, the following chain of equalities holds
\allowdisplaybreaks
\begin{align*}
f^{[ 1 ] @}_{s} \left(
\mathrm{ip}^{([ 1 ], 0)\sharp}_{\boldsymbol{\mathcal{A}}^{(1)}, s}\left(
P
\right)
\right)
&=
f^{[ 1 ] @}_{s} \left(
\left[
\mathrm{ip}^{(1, 0)\sharp}_{\boldsymbol{\mathcal{A}}^{(1)}, s}\left(
P
\right)
\right]_{s}
\right)
\tag{1}
\\
&=
\left[
f^{(1)\flat}_{s} \left(
\mathrm{ip}^{(1, 0)\sharp}_{\boldsymbol{\mathcal{A}}^{(1)}, s}\left(
P
\right)
\right)
\right]_{\varphi(s)}
\tag{2}
\\
&=
\left[
\mathrm{ip}^{(1, 0)\sharp}_{\boldsymbol{\mathcal{B}}^{(1)}, \varphi(s)}\left(
f^{(0)\sharp}_{s} \left(
P
\right)
\right)
\right]_{\varphi(s)}
\tag{3}
\\
&=
\mathrm{ip}^{([ 1 ], 0)\sharp}_{\boldsymbol{\mathcal{B}}^{(1)}, \varphi(s)}\left(
f^{(0)\sharp}_{s} \left(
P
\right)
\right)
\tag{4}
\end{align*}

The first equality unravels the definition of the mapping $\mathrm{ip}^{([1], 0)}_{\boldsymbol{\mathcal{A}}^{(1)}}$, introduced in Definition~\ref{DCHUZ};
the second equality unravels the definition of the first-order quotient path extension mapping $f^{[ 1 ] @}$ at a path class, introduced in Definition~\ref{DQPthExt};
the third equality follows from the definition of the path extension mapping, introduced in Proposition~\ref{PPthExt};
finally, the last equality recovers the definition of the mapping $\mathrm{ip}^{([1], 0)}_{\boldsymbol{\mathcal{B}}^{(1)}, \varphi}$, introduced in Definition~\ref{DCHUZ}.
\end{proof}

\begin{proposition}
\label{PQPthExtScTg}
Let $\mathbf{f}^{(1)}=(\varphi, c, (f^{(i)})_{i\in 2})$ be a morphism of first-order many-sorted rewriting systems from $\boldsymbol{\mathcal{A}}^{(1)}$ to $\boldsymbol{\mathcal{B}}^{(1)}$. Then the following equalities holds
\begin{align*}
\mathrm{sc}_{\boldsymbol{\mathcal{B}}^{(1)},\varphi}^{(0,[1])} \circ f^{[1] @}
&=
f^{(0)\sharp}\circ\mathrm{sc}_{\boldsymbol{\mathcal{A}}^{(1)}}^{(0,[1])}
&&\mbox{and}&
\mathrm{tg}_{\boldsymbol{\mathcal{B}}^{(1)},\varphi}^{(0,[1])} \circ f^{[1] @}
&=
f^{(0)\sharp}\circ\mathrm{tg}_{\boldsymbol{\mathcal{A}}^{(1)}}^{(0,[1])}.
\end{align*}
\end{proposition}

\begin{proof}
For every sort $s$ in $S$ and every path class $[ \mathfrak{P} ]_{s}$ in $[ \mathrm{Pth}_{\boldsymbol{\mathcal{A}}^{(1)}} ]_{s}$, the following chain of equalities holds
\allowdisplaybreaks
\begin{align*}
\mathrm{sc}^{(0,[1])}_{\boldsymbol{\mathcal{B}}^{(1)}, \varphi(s)} \left(
f^{[ 1 ] @}_{s} \left(
\left[ \mathfrak{P} \right]_{s}
\right)
\right)
&=
\mathrm{sc}^{(0,[1])}_{\boldsymbol{\mathcal{B}}^{(1)}, \varphi(s)} \left(
\left[
f^{(1)\flat}_{s} \left(
\mathfrak{P}
\right)
\right]_{\varphi(s)}
\right)
\tag{1}
\\
&=
\mathrm{sc}^{(0,1)}_{\boldsymbol{\mathcal{B}}^{(1)}, \varphi(s)} \left(
f^{(1)\flat}_{s} \left(
\mathfrak{P}
\right)
\right)
\tag{2}
\\
&=
f^{(0)\sharp}_{s} \left(
\mathrm{sc}^{(0,1)}_{\boldsymbol{\mathcal{A}}^{(1)}, s} \left(
\mathfrak{P}
\right)
\right)
\tag{3}
\\
&=
f^{(0)\sharp}_{s} \left(
\mathrm{sc}^{(0,[1])}_{\boldsymbol{\mathcal{A}}^{(1)}, s} \left(
\left[
\mathfrak{P}
\right]_{s}
\right)
\right)
\tag{4}
\end{align*}

The first equality unravels the definition of the first-order quotient path extension mapping $f^{[ 1 ] @}$ at a path class, introduced in Definition~\ref{DQPthExt};
the second equality unravels the definition of the mapping $\mathrm{sc}^{(0,[1])}_{\boldsymbol{\mathcal{B}}^{(1)}, \varphi}$, introduced in Definition~\ref{DCHUZ};
the third equality follows from Proposition~\ref{PPthExt};
finally, the last equality recovers the definition of the mapping $\mathrm{sc}^{(0,[1])}_{\boldsymbol{\mathcal{A}}^{(1)}}$, introduced in Definition~\ref{DCHUZ}.

A similar argument applies to the $(0,[ 1 ])$-target mapping.
\end{proof}

\begin{proposition}
\label{PQPthExtIp}
Let $\mathbf{f}^{(1)}=(\varphi, c, (f^{(i)})_{i\in 2})$ be a morphism of first-order many-sorted rewriting systems from $\boldsymbol{\mathcal{A}}^{(1)}$ to $\boldsymbol{\mathcal{B}}^{(1)}$. Then the following equality holds
$$
f^{[1] @} \circ \mathrm{ip}_{\boldsymbol{\mathcal{A}}^{(1)}}^{([1], 0)\sharp}
=
\mathrm{ip}_{\boldsymbol{\mathcal{B}}^{(1)}, \varphi}^{([1], 0)\sharp} \circ f^{(0)\sharp}.
$$
\end{proposition}

\begin{proof}
For every sort $s$ in $S$ and every term $P$ in $\T_{\Sigma}(X)$, the following chain of equalities holds
\allowdisplaybreaks
\begin{align*}
f^{[ 1 ] @}_{s} \left(
\mathrm{ip}^{([ 1 ], 0)\sharp}_{\boldsymbol{\mathcal{A}}^{(1)}, s}\left(
P
\right)
\right)
&=
f^{[ 1 ] @}_{s} \left(
\left[
\mathrm{ip}^{(1, 0)\sharp}_{\boldsymbol{\mathcal{A}}^{(1)}, s}\left(
P
\right)
\right]_{s}
\right)
\tag{1}
\\
&=
\left[
f^{(1)\flat}_{s} \left(
\mathrm{ip}^{(1, 0)\sharp}_{\boldsymbol{\mathcal{A}}^{(1)}, s}\left(
P
\right)
\right)
\right]_{\varphi(s)}
\tag{2}
\\
&=
\left[
\mathrm{ip}^{(1, 0)\sharp}_{\boldsymbol{\mathcal{B}}^{(1)}, \varphi(s)}\left(
f^{[1]@}_{s} \left(
P
\right)
\right)
\right]_{\varphi(s)}
\tag{3}
\\
&=
\mathrm{ip}^{([ 1 ], 0)\sharp}_{\boldsymbol{\mathcal{B}}^{(1)}, \varphi(s)}\left(
f^{[1]@}_{s} \left(
P
\right)
\right)
\tag{4}
\end{align*}

The first equality unravels the definition of the mapping $\mathrm{ip}^{([1], 0)}_{\boldsymbol{\mathcal{A}}^{(1)}}$, introduced in Definition~\ref{DCHUZ};
the second equality unravels the definition of the first-order quotient path extension mapping $f^{[ 1 ] @}$ at a path class, introduced in Definition~\ref{DQPthExt};
the third equality follows from Proposition~\ref{PPthExt};
finally, the last equality recovers the definition of the mapping $\mathrm{ip}^{([1], 0)}_{\boldsymbol{\mathcal{B}}^{(1)}, \varphi}$, introduced in Definition~\ref{DCHUZ}.
\end{proof}			
\chapter{The category of first-order many-sorted rewriting systems}\label{S3D}

In this chapter we introduce the notions of first-order identity morphism at a first-order many-sorted rewriting system and that of first-order composition morphism of first-order morphisms. We then define the first-order equivalence relation between first-order morphisms to show that first-order many-sorted rewriting systems with equivalence classes of first-order morphisms form a category denoted $\mathsf{Rws}_{\mathfrak{d}}^{[1]}$. Moreover, we introduce the notion of first-order tower associated to a first-order rewriting system and that of morphism of first-order towers. Then showing that the first-order towers and morphisms between them form a category denoted $\mathsf{Tw}_{\mathfrak{d}}^{[1]}$. Indeed we show that the categories $\mathsf{Rws}_{\mathfrak{d}}^{[1]}$ and $\mathsf{Tw}_{\mathfrak{d}}^{[1]}$ are isomorphic. 

\section{The category $\mathsf{Rws}_{\mathfrak{d}}^{[1]}$}

We begin by defining the identity morphism at a first-order rewriting system and studying its properties.

\begin{definition}
\label{DIdRws1}
Let $\boldsymbol{\mathcal{A}}^{(1)} = (\boldsymbol{\mathcal{A}}^{(0)}, \mathcal{A}^{(1)})$ be a first-order many-sorted rewriting system. The \emph{first-order identity morphism} at $\boldsymbol{\mathcal{A}}^{(1)}$, denoted by $\mathrm{id}^{\boldsymbol{\mathcal{A}}^{(1)}}$, is given by 
$$
\mathrm{id}^{\boldsymbol{\mathcal{A}}^{(1)}}
=
\left(
\mathrm{id}^{\boldsymbol{\mathcal{A}}^{(0)}}, 
\mathrm{ech}^{(1,\mathcal{A}^{(1)})}_{\boldsymbol{\mathcal{A}}^{(1)}}
\right).
$$
where $\mathrm{id}^{\boldsymbol{\mathcal{A}}^{(0)}}$ is the identity zeroth-order morphism at the zeroth-order many-sorted rewriting system $\boldsymbol{\mathcal{A}}^{(0)}$, introduced in Definition~\ref{DIdRws0}, and $\mathrm{ech}^{(1,\mathcal{A}^{(1)})}_{\boldsymbol{\mathcal{A}}^{(1)}}$ is the $S$-sorted first-order echelon mapping associated with the many-sorted rewriting system $\boldsymbol{\mathcal{A}}^{(1)}$, introduced in Definition~\ref{DEch}.
\end{definition}

\begin{proposition}
\label{PPthExtEch}
Let $\boldsymbol{\mathcal{A}}^{(1)}$ be a first-order rewriting system and let $\mathrm{id}^{\boldsymbol{\mathcal{A}}^{(1)}}$ be its first-order identity morphism. Thus, $\mathrm{ech}^{(1,\mathcal{A}^{(1)})\flat}_{\boldsymbol{\mathcal{A}}^{(1)}} = \mathrm{ip}_{\boldsymbol{\mathcal{A}}^{(1)}}^{(1,X)@} \circ \mathrm{CH}^{(1)}_{\boldsymbol{\mathcal{A}}^{(1)}}$.
\end{proposition}

\begin{proof}
We want to show that, for every sort $s \in S$ and every path $\mathfrak{P}$ in $\mathrm{Pth}_{\boldsymbol{\mathcal{A}}^{(1)},s}$, the path $\mathrm{ech}^{(1,\mathcal{A}^{(1)})\flat}_{\boldsymbol{\mathcal{A}}^{(1)}} (\mathfrak{P})$ is equal to $\mathrm{ip}_{\boldsymbol{\mathcal{A}}^{(1)}}^{(1,X)@} \circ \mathrm{CH}^{(1)}_{\boldsymbol{\mathcal{A}}^{(1)}} (\mathfrak{P})$.

We prove the statement by Artinian induction on $(\coprod\mathrm{Pth}_{\boldsymbol{\mathcal{A}}^{(1)}}, \leq_{\mathbf{Pth}_{\boldsymbol{\mathcal{A}}^{(1)}}})$.

\textsf{Base step of the Artinian induction.}

Let $(\mathfrak{P}, s)$ be a minimal element in $(\coprod\mathrm{Pth}_{\boldsymbol{\mathcal{A}}^{(1)}}, \leq_{\mathbf{Pth}_{\boldsymbol{\mathcal{A}}^{(1)}}})$. Then by Proposition~\ref{PMinimal}, it follows that $\mathfrak{P}$ is either (1) an $(1,0)$-identity path or (2) an echelon.

If~(1), i.e., if $\mathfrak{P}$ is an $(1,0)$-identity path, then $\mathfrak{P} = \mathrm{ip}^{(1,0)\sharp}_{\boldsymbol{\mathcal{A}}^{(1)}, s}(P)$ for some term $P \in \mathrm{T}_{\Sigma}(X)_{s}$. Thus, the following chain of equalities holds
\allowdisplaybreaks
\begin{align*}
\mathrm{ech}^{(1,\mathcal{A}^{(1)})\flat}_{\boldsymbol{\mathcal{A}}^{(1)}, s}(\mathfrak{P})
&=
\mathrm{ech}^{(1,\mathcal{A}^{(1)})\flat}_{\boldsymbol{\mathcal{A}}^{(1)}, s}\left(
\mathrm{ip}^{(1,0)\sharp}_{\boldsymbol{\mathcal{A}}^{(1)}, s}\left(
P
\right)
\right)
\tag{1}
\\
&=
\mathrm{ip}^{(1,0)\sharp}_{\boldsymbol{\mathcal{A}}^{(1)}, s}\left(
\eta^{X\sharp}\left(
P
\right)
\right)
\tag{2}
\\
&=
\mathrm{ip}^{(1,0)\sharp}_{\boldsymbol{\mathcal{A}}^{(1)}, s}\left(
P
\right)
\tag{3}
\\
&=
\eta^{(1,\mathbf{Pth}_{\boldsymbol{\mathcal{A}}^{(1)}})}\left(
\mathrm{ip}^{(1,0)\sharp}_{\boldsymbol{\mathcal{A}}^{(1)}, s}\left(
P
\right)
\right)
\tag{4}
\\
&=
\mathrm{ip}^{(1,X)@}_{\boldsymbol{\mathcal{A}}^{(1)}, s}\left(
\mathrm{CH}^{(1)}_{\boldsymbol{\mathcal{A}}^{(1)}, s}\left(
\mathrm{ip}^{(1,0)\sharp}_{\boldsymbol{\mathcal{A}}^{(1)}, s}\left(
P
\right)
\right)
\right)
\tag{5}
\\
&=
\mathrm{ip}^{(1,X)@}_{\boldsymbol{\mathcal{A}}^{(1)}, s}\left(
\mathrm{CH}^{(1)}_{\boldsymbol{\mathcal{A}}^{(1)}, s}\left(
\mathfrak{P}
\right)
\right).
\tag{6}
\end{align*}

In the just stated chain of equalities, the first equality unravels the definition of the path $\mathfrak{P}$;
the second equality unravels the definition of the path extension mapping at a $(1,0)$-identity path, introduced in Proposition~\ref{PPthExt};
the third equality follows from the fact that $\eta^{X\sharp}$ is equal to $\mathrm{id}^{\mathrm{T}_{\Sigma}(X)}$;
the fourth equality follows from the fact that we identify $\eta^{(1,\mathbf{Pth}_{\boldsymbol{\mathcal{A}}^{(1)}})}(\mathfrak{P})$ with $\mathfrak{P}$;
the fifth equality follows from Proposition~\ref{PIpUZ};
finally, the last equality recovers the definition of the path $\mathfrak{P}$.

If~(2), i.e., if $\mathfrak{P}$ is an echelon associated to a rewrite rule $\mathfrak{p} = (M,N)$ in $\mathcal{A}^{(1)}_{s}$, that is, if $\mathfrak{P}$ has the form 
$$
\xymatrix@C=55pt{
\mathfrak{P}: M
\ar[r]^-{\text{\Small{$(\mathfrak{p},\mathrm{id}^{\mathrm{T}_{\Sigma}(X)_{s}})$}}}
&
N
},
$$
then the following chain of equalities holds
\begin{align*}
\mathrm{ech}^{(1,\mathcal{A}^{(1)})\flat}_{\boldsymbol{\mathcal{A}}^{(1)}, s}(\mathfrak{P})
&=
\mathrm{ech}^{(1,\mathcal{A}^{(1)})\flat}_{\boldsymbol{\mathcal{A}}^{(1)}, s}\left(
\mathrm{ech}^{(1,\mathcal{A}^{(1)})}_{\boldsymbol{\mathcal{A}}^{(1)}, s}\left(\mathfrak{p}\right)
\right)
\tag{1}
\\
&=
\mathrm{ech}^{(1,\mathcal{A}^{(1)})}_{\boldsymbol{\mathcal{A}}^{(1)}, s}\left(
\mathfrak{p}
\right)
\tag{2}
\\
&=
\mathfrak{p}^{\mathbf{Pth}_{\boldsymbol{\mathcal{A}}^{(1)}}}
\tag{3}
\\
&=
\mathfrak{p}^{\mathbf{F}_{\Sigma^{\boldsymbol{\mathcal{A}}^{(1)}}}(\mathbf{Pth}_{\boldsymbol{\mathcal{A}}^{(1)}})}
\tag{4}
\\
&=
\mathrm{ip}^{(1,X)@}_{\boldsymbol{\mathcal{A}}^{(1)}, s}\left(
\mathfrak{p}^{\mathbf{T}_{\Sigma^{\boldsymbol{\mathcal{A}}^{(1)}}}(X)}
\right)
\tag{5}
\\
&=
\mathrm{ip}^{(1,X)@}_{\boldsymbol{\mathcal{A}}^{(1)}, s}\left(
\mathrm{CH}^{(1)}_{\boldsymbol{\mathcal{A}}^{(1)}, s} \left(
\mathrm{ech}^{(1,\mathcal{A}^{(1)})}_{\boldsymbol{\mathcal{A}}^{(1)}, s}\left(\mathfrak{p}\right)
\right)
\right)
\tag{6}
\\
&=
\mathrm{ip}^{(1,X)@}_{\boldsymbol{\mathcal{A}}^{(1)}, s}\left(
\mathrm{CH}^{(1)}_{\boldsymbol{\mathcal{A}}^{(1)}, s} \left(
\mathfrak{P}
\right)
\right).
\tag{7}
\end{align*}

In the just stated chain of equalities, the first equality unravels the definition of $\mathfrak{P}$;
the second equality unravels the definition of the path extension mapping at an echelon, introduced in Proposition~\ref{PPthExt};
the third equality recovers the interpretation of the constant symbol $\mathfrak{p}$ in the many-sorted partial $\Sigma^{\boldsymbol{\mathcal{A}}^{(1)}}$-algebra introduced in Proposition~\ref{PPthCatAlg};
the fourth equality holds since, by Proposition~\ref{PPthCatAlg}, we have an interpretation of the constant symbol $\mathfrak{p}$;
the fifth equality holds because, according to Definition~\ref{DIp}, $\mathrm{ip}^{(1,X)@}_{\boldsymbol{\mathcal{A}}^{(1)}}$ is a many-sorted $\Sigma^{\boldsymbol{\mathcal{A}}^{(1)}}$-homomorphism;
the sixth equality recovers the definition of the Curry-Howard mapping at an echelon, introduced in Definition~\ref{DCH};
finally, the last equality recovers the definition of $\mathfrak{P}$.

This proves the base case of the induction.

\textsf{Inductive step of the Artinian induction.}

Let $(\mathfrak{P},s)$ be a non-minimal element of $(\coprod\mathrm{Pth}_{\boldsymbol{\mathcal{A}}^{(1)}}, \leq_{\mathbf{Pth}_{\boldsymbol{\mathcal{A}}^{(1)}}})$. We can assume that $\mathfrak{P}$ is not a $(1,0)$-identity path, since that case has already been considered. Let us suppose that, for every sort $t \in S$ and every path $\mathfrak{Q} \in \mathrm{Pth}_{\boldsymbol{\mathcal{A}}^{(1)}, t}$, if $(\mathfrak{Q}, t) <_{\mathbf{Pth}_{\boldsymbol{\mathcal{A}}^{(1)}}} (\mathfrak{P}, s)$, then the value of $\mathrm{ech}^{(1,\mathcal{A}^{(1)})\flat}_{\boldsymbol{\mathcal{A}}^{(1)}, t}(\mathfrak{Q})$ is equal to $\mathrm{ip}^{(1,X)@}_{\boldsymbol{\mathcal{A}}^{(1)}, t}(
\mathrm{CH}^{(1)}_{\boldsymbol{\mathcal{A}}^{(1)}, t}(
\mathfrak{Q}
)
)$. 

By Lemma~\ref{LOrdI}, we have that $\mathfrak{P}$ is either~(1) a path of length strictly greater that one containing at least one echelon or~(2) an echelonless path.

If~(1), i.e., if $\mathfrak{P}$ is a path of legth strictly greater than one containing at least one echelon. Then, we let $i \in \bb{\mathfrak{P}}$ be the first index for which the one-step subpath $\mathfrak{P}^{i,i}$ is an echelon. We distinguish two cases accordingly, either~(1.1) $i=0$ or~(1.2) $i>0$. 

If~(1.1), i.e., if $\mathfrak{P}$ is a path of length strictly greater having its first echelon on its first step, then the following chain of equalities holds
\begin{flushleft}
$
\mathrm{ech}^{(1,\mathcal{A}^{(1)})\flat}_{\boldsymbol{\mathcal{A}}^{(1)}, s}(\mathfrak{P})
$
\allowdisplaybreaks
\begin{align*}
&=
\mathrm{ech}^{(1,\mathcal{A}^{(1)})\flat}_{\boldsymbol{\mathcal{A}}^{(1)}, s}\left(
\mathfrak{P}^{1, \bb{\mathfrak{P}}-1}
\right)
\circ^{0\mathbf{Pth}_{\boldsymbol{\mathcal{A}}^{(1)}}}_{s}
\mathrm{ech}^{(1,\mathcal{A}^{(1)})\flat}_{\boldsymbol{\mathcal{A}}^{(1)}, s}\left(
\mathfrak{P}^{0,0}
\right)
\tag{1}
\\
&=
\mathrm{ip}^{(1,X)@}_{\boldsymbol{\mathcal{A}}^{(1)}, s}\left(
\mathrm{CH}^{(1)}_{\boldsymbol{\mathcal{A}}^{(1)}, s} \left(
\mathfrak{P}^{1, \bb{\mathfrak{P}}-1}
\right)
\right)
\circ^{0\mathbf{Pth}_{\boldsymbol{\mathcal{A}}^{(1)}}}_{s}
\mathrm{ip}^{(1,X)@}_{\boldsymbol{\mathcal{A}}^{(1)}, s}\left(
\mathrm{CH}^{(1)}_{\boldsymbol{\mathcal{A}}^{(1)}, s} \left(
\mathfrak{P}^{0,0}
\right)
\right)
\tag{2}
\\
&=
\resizebox{.89\textwidth}{!}{$
\mathrm{ip}^{(1,X)@}_{\boldsymbol{\mathcal{A}}^{(1)}, s}\left(
\mathrm{CH}^{(1)}_{\boldsymbol{\mathcal{A}}^{(1)}, s} \left(
\mathfrak{P}^{1, \bb{\mathfrak{P}}-1}
\right)
\right)
\circ^{0\mathbf{F}_{\Sigma^{\boldsymbol{\mathcal{A}}^{(1)}}}(\mathbf{Pth}_{\boldsymbol{\mathcal{A}}^{(1)}})}_{s}
\mathrm{ip}^{(1,X)@}_{\boldsymbol{\mathcal{A}}^{(1)}, s}\left(
\mathrm{CH}^{(1)}_{\boldsymbol{\mathcal{A}}^{(1)}, s} \left(
\mathfrak{P}^{0,0}
\right)
\right)
$}
\tag{3}
\\
&=
\mathrm{ip}^{(1,X)@}_{\boldsymbol{\mathcal{A}}^{(1)}, s}\left(
\mathrm{CH}^{(1)}_{\boldsymbol{\mathcal{A}}^{(1)}, s} \left(
\mathfrak{P}^{1, \bb{\mathfrak{P}}-1}
\right)
\circ^{0\mathbf{T}_{\Sigma^{\boldsymbol{\mathcal{A}}^{(1)}}}(X)}_{s}
\mathrm{CH}^{(1)}_{\boldsymbol{\mathcal{A}}^{(1)}, s} \left(
\mathfrak{P}^{0,0}
\right)
\right)
\tag{4}
\\
&=
\mathrm{ip}^{(1,X)@}_{\boldsymbol{\mathcal{A}}^{(1)}, s}\left(
\mathrm{CH}^{(1)}_{\boldsymbol{\mathcal{A}}^{(1)}, s} \left(
\mathfrak{P}
\right)
\right).
\tag{5}
\end{align*}
\end{flushleft}

In the just stated chain of equalities, the first equality unravels the definition of the path extension mapping at a path of length strictly greater than one with an echelon on its first step, introduced in Proposition~\ref{PPthExt};
the second equality follows taking into account that $(\mathfrak{P}^{(1,\bb{\mathfrak{P}}-1}, s)$ and $(\mathfrak{P}^{0,0}, s)$ are strictly smaller that $(\mathfrak{P}, s)$ with respect to $\prec_{\mathbf{Pth}_{\boldsymbol{\mathcal{A}}^{(1)}}}$;
the third equality holds because, according to Proposition~\ref{PIpCH}, we have that 
$\mathrm{ip}^{(1,X)@}_{\boldsymbol{\mathcal{A}}^{(1)}, s}(
\mathrm{CH}^{(1)}_{\boldsymbol{\mathcal{A}}^{(1)}, s} (
\mathfrak{P}^{1, \bb{\mathfrak{P}}-1}
)
)$
and 
$
\mathrm{ip}^{(1,X)@}_{\boldsymbol{\mathcal{A}}^{(1)}, s}(
\mathrm{CH}^{(1)}_{\boldsymbol{\mathcal{A}}^{(1)}, s}(
\mathfrak{P}^{0,0}
)
)
$
are paths in $[\mathfrak{P}^{1, \bb{\mathfrak{P}}-1}]_{s}$ and $[\mathfrak{P}^{0,0}]_{s}$, respectively. Hence the interpretation of the $0$-composition operation symbol $\circ^{0}_{s}$ in $\mathbf{F}_{\Sigma^{\boldsymbol{\mathcal{A}}^{(1)}}}(\mathbf{Pth}_{\boldsymbol{\mathcal{A}}^{(1)}})$ becomes that of $\mathbf{Pht}_{\boldsymbol{\mathcal{A}}^{(1)}}$.
The fourth equality holds since, according to Definition~\ref{DIp}, $\mathrm{ip}^{(1,X)@}_{\boldsymbol{\mathcal{A}}^{(1)}}$ is a many-sorted $\Sigma^{\boldsymbol{\mathcal{A}}^{(1)}}$-homomorphism;
finally, the last equality recovers the definition of the Curry-Howard mapping at a path of length strictly greater than one with an echelon on its first step, introduced in Definition~\ref{DCH}.

If~(1.2), i.e., if $\mathfrak{P}$ is a path of length strictly greater having its first echelon at position $i \in \bb{\mathfrak{P}}$, with $i > 0$, then the following chain of equalities holds
\begin{flushleft}
$
\mathrm{ech}^{(1,\mathcal{A}^{(1)})\flat}_{\boldsymbol{\mathcal{A}}^{(1)}, s}(\mathfrak{P})
$
\allowdisplaybreaks
\begin{align*}
&=
\mathrm{ech}^{(1,\mathcal{A}^{(1)})\flat}_{\boldsymbol{\mathcal{A}}^{(1)}, s}\left(
\mathfrak{P}^{i, \bb{\mathfrak{P}}-1}
\right)
\circ^{0\mathbf{Pth}_{\boldsymbol{\mathcal{B}}^{(1)}}}_{\varphi(s)}
\mathrm{ech}^{(1,\mathcal{A}^{(1)})\flat}_{\boldsymbol{\mathcal{A}}^{(1)}, s}\left(
\mathfrak{P}^{0,i-1}
\right)
\tag{1}
\\
&=
\mathrm{ip}^{(1,X)@}_{\boldsymbol{\mathcal{A}}^{(1)}, s}\left(
\mathrm{CH}^{(1)}_{\boldsymbol{\mathcal{A}}^{(1)}, s} \left(
\mathfrak{P}^{i, \bb{\mathfrak{P}}-1}
\right)
\right)
\circ^{0\mathbf{Pth}_{\boldsymbol{\mathcal{A}}^{(1)}}}_{s}
\mathrm{ip}^{(1,X)@}_{\boldsymbol{\mathcal{A}}^{(1)}, s}\left(
\mathrm{CH}^{(1)}_{\boldsymbol{\mathcal{A}}^{(1)}, s} \left(
\mathfrak{P}^{0,i-1}
\right)
\right)
\tag{2}
\\
&=
\resizebox{.89\textwidth}{!}{$
\mathrm{ip}^{(1,X)@}_{\boldsymbol{\mathcal{A}}^{(1)}, s}\left(
\mathrm{CH}^{(1)}_{\boldsymbol{\mathcal{A}}^{(1)}, s} \left(
\mathfrak{P}^{i, \bb{\mathfrak{P}}-1}
\right)
\right)
\circ^{0\mathbf{F}_{\Sigma^{\boldsymbol{\mathcal{A}}^{(1)}}}(\mathbf{Pth}_{\boldsymbol{\mathcal{A}}^{(1)}})}_{s}
\mathrm{ip}^{(1,X)@}_{\boldsymbol{\mathcal{A}}^{(1)}, s}\left(
\mathrm{CH}^{(1)}_{\boldsymbol{\mathcal{A}}^{(1)}, s} \left(
\mathfrak{P}^{0,i-1}
\right)
\right)
$}
\tag{3}
\\
&=
\mathrm{ip}^{(1,X)@}_{\boldsymbol{\mathcal{A}}^{(1)}, s}\left(
\mathrm{CH}^{(1)}_{\boldsymbol{\mathcal{A}}^{(1)}, s} \left(
\mathfrak{P}^{i, \bb{\mathfrak{P}}-1}
\right)
\circ^{0\mathbf{T}_{\Sigma^{\boldsymbol{\mathcal{A}}^{(1)}}}(X)}_{s}
\mathrm{CH}^{(1)}_{\boldsymbol{\mathcal{A}}^{(1)}, s} \left(
\mathfrak{P}^{0,i-1}
\right)
\right)
\tag{4}
\\
&=
\mathrm{ip}^{(1,X)@}_{\boldsymbol{\mathcal{A}}^{(1)}, s}\left(
\mathrm{CH}^{(1)}_{\boldsymbol{\mathcal{A}}^{(1)}, s} \left(
\mathfrak{P}
\right)
\right).
\tag{5}
\end{align*}
\end{flushleft}

In the just stated chain of equalities, the first equality unravels the definition of the path extension mapping at a path of length strictly greater than one with an echelon at position $i$, introduced in Proposition~\ref{PPthExt};
the second equality follows taking into account that $(\mathfrak{P}^{(i,\bb{\mathfrak{P}}-1}, s)$ and $(\mathfrak{P}^{0,i-1}, s)$ are strictly smaller that $(\mathfrak{P}, s)$ with respect to $\prec_{\mathbf{Pth}_{\boldsymbol{\mathcal{A}}^{(1)}}}$;
the third equality holds because, according to Proposition~\ref{PIpCH}, we have that 
$\mathrm{ip}^{(1,X)@}_{\boldsymbol{\mathcal{A}}^{(1)}, s}(
\mathrm{CH}^{(1)}_{\boldsymbol{\mathcal{A}}^{(1)}, s} (
\mathfrak{P}^{i, \bb{\mathfrak{P}}-1}
)
)$
and 
$
\mathrm{ip}^{(1,X)@}_{\boldsymbol{\mathcal{A}}^{(1)}, s}(
\mathrm{CH}^{(1)}_{\boldsymbol{\mathcal{A}}^{(1)}, s}(
\mathfrak{P}^{0,i-1}
)
)
$
are paths in $[\mathfrak{P}^{i, \bb{\mathfrak{P}}-1}]_{s}$ and $[\mathfrak{P}^{0,i-1}]_{s}$, respectively. Hence the interpretation of the $0$-composition operation symbol $\circ^{0}_{s}$ in $\mathbf{F}_{\Sigma^{\boldsymbol{\mathcal{A}}^{(1)}}}(\mathbf{Pth}_{\boldsymbol{\mathcal{A}}^{(1)}})$ becomes that of $\mathbf{Pth}_{\boldsymbol{\mathcal{A}}^{(1)}}$.
The fourth equality holds since, according to Definition~\ref{DIp}, $\mathrm{ip}^{(1,X)@}_{\boldsymbol{\mathcal{A}}^{(1)}}$ is a many-sorted $\Sigma^{\boldsymbol{\mathcal{A}}^{(1)}}$-homomorphism;
finally, the last equality recovers the definition of the Curry-Howard mapping at a path of length strictly greater than one with an echelon at position $i$, introduced in Definition~\ref{DCH}.

If~(2), i.e., if $\mathfrak{P}$ is an echelonless path then, according to Lemma~\ref{LPthHeadCt}, there exsits a unique word $s \in S^{\star}-\{\lambda\}$ and a unique operation symbol $\sigma$ in $\Sigma_{\mathbf{s}, s}$ associated to $\mathfrak{P}$. Let $(\mathfrak{P}_{j})_{j \in \bb{\mathbf{s}}}$ be the family of paths we can extract from $\mathfrak{P}$ in virtue of Lemma~\ref{LPthExtract}. Then, the following chain of equalities holds
\begin{align*}
\mathrm{ech}^{(1,\mathcal{A}^{(1)})\flat}_{\boldsymbol{\mathcal{A}}^{(1)}, s}(\mathfrak{P})
&=
\sigma^{\mathbf{Pth}_{\boldsymbol{\mathcal{A}}^{(1)}}^{\mathrm{id}^{\boldsymbol{\mathcal{A}}^{(1)}}(0,1)}}\left(
\left(
\mathrm{ech}^{(1,\mathcal{A}^{(1)})\flat}_{\boldsymbol{\mathcal{A}}^{(1)}, s_{j}}(\mathfrak{P}_{j})
\right)_{j \in \bb{\mathbf{s}}}
\right)
\tag{1}
\\
&=
\sigma^{\mathbf{Pth}_{\boldsymbol{\mathcal{A}}^{(1)}}^{\mathrm{id}^{\boldsymbol{\mathcal{A}}^{(1)}}(0,1)}}\left(
\left(
\mathrm{ip}^{(1,X)@}_{\boldsymbol{\mathcal{A}}^{(1)}, s}\left(
\mathrm{CH}^{(1)}_{\boldsymbol{\mathcal{A}}^{(1)}, s} \left(
\mathfrak{P}_{j}
\right)
\right)
\right)_{j \in \bb{\mathbf{s}}}
\right)
\tag{2}
\\
&=
\sigma^{\mathbf{Pth}_{\boldsymbol{\mathcal{A}}^{(1)}}^{(0,1)}}\left(
\left(
\mathrm{ip}^{(1,X)@}_{\boldsymbol{\mathcal{A}}^{(1)}, s}\left(
\mathrm{CH}^{(1)}_{\boldsymbol{\mathcal{A}}^{(1)}, s} \left(
\mathfrak{P}_{j}
\right)
\right)
\right)_{j \in \bb{\mathbf{s}}}
\right)
\tag{3}
\\
&=
\sigma^{\mathbf{F}_{\Sigma^{\boldsymbol{\mathcal{A}}^{(1)}}}(\mathbf{Pth}_{\boldsymbol{\mathcal{A}}^{(1)}})}\left(
\left(
\mathrm{ip}^{(1,X)@}_{\boldsymbol{\mathcal{A}}^{(1)}, s}\left(
\mathrm{CH}^{(1)}_{\boldsymbol{\mathcal{A}}^{(1)}, s} \left(
\mathfrak{P}_{j}
\right)
\right)
\right)_{j \in \bb{\mathbf{s}}}
\right)
\tag{4}
\\
&=
\mathrm{ip}^{(1,X)@}_{\boldsymbol{\mathcal{A}}^{(1)}, s}\left(
\sigma^{\mathbf{T}_{\Sigma^{\boldsymbol{\mathcal{A}}^{(1)}}}(X)}\left(
\left(
\mathrm{CH}^{(1)}_{\boldsymbol{\mathcal{A}}^{(1)}, s} \left(
\mathfrak{P}_{j}
\right)
\right)
\right)_{j \in \bb{\mathbf{s}}}
\right)
\tag{5}
\\
&=
\mathrm{ip}^{(1,X)@}_{\boldsymbol{\mathcal{A}}^{(1)}, s}\left(
\mathrm{CH}^{(1)}_{\boldsymbol{\mathcal{A}}^{(1)}, s} \left(
\mathfrak{P}
\right)
\right).
\tag{6}
\end{align*}

In the just stated chain of equalities, the first equality unravels the definition of the path extension mapping at an echelonless path, introduced in Proposition~\ref{PPthExt};
the second equality follows taking into account that, for every $j\in\bb{\mathbf{s}}$, $(\mathfrak{P}_{j}, s_{j})$ is strictly smaller than $(\mathfrak{P}, s)$ with respect to $\prec_{\mathbf{Pth}_{\boldsymbol{\mathcal{A}}^{(1)}}}$;
the third equality follows from the fact that, according to Proposition~\ref{PAlgFun}, $\mathrm{Alg}_{\mathfrak{d}}$ is a functor, thus $\mathbf{Pth}_{\boldsymbol{\mathcal{A}}^{(1)}}^{\mathrm{id}^{\boldsymbol{\mathcal{A}}^{(1)}}(0,1)} 
= 
(\mathrm{id}^{\boldsymbol{\Sigma}})^{*}_{\mathfrak{d}}(\mathbf{Pth}_{\boldsymbol{\mathcal{A}}^{(1)}}^{(0,1)})
=
\mathbf{Pth}_{\boldsymbol{\mathcal{A}}^{(1)}}^{(0,1)}$;
the fourth equality holds because, according to Proposition~\ref{PIpCH}, we have that, for every $j \in \bb{\mathbf{s}}$,
$\mathrm{ip}^{(1,X)@}_{\boldsymbol{\mathcal{A}}^{(1)}, s}(
\mathrm{CH}^{(1)}_{\boldsymbol{\mathcal{A}}^{(1)}, s} (
\mathfrak{P}_{j}
)
)$
is a path in $[\mathfrak{P}_{j}]_{s_{j}}$. Hence the interpretation of the operation symbol $\sigma$ in $\mathbf{F}_{\Sigma^{\boldsymbol{\mathcal{A}}^{(1)}}}(\mathbf{Pth}_{\boldsymbol{\mathcal{A}}^{(1)}})$ becomes that of $\mathbf{Pht}_{\boldsymbol{\mathcal{A}}^{(1)}}$
the fifth equality holds since, according to Definition~\ref{DIp}, $\mathrm{ip}^{(1,X)@}_{\boldsymbol{\mathcal{A}}^{(1)}}$ is a many-sorted $\Sigma^{\boldsymbol{\mathcal{A}}^{(1)}}$-homomorphism;
finally, the last equality recovers the definition of the Curry-Howard mapping at an echelonless path, introduced in Definition~\ref{DCH}.

This proves the inductive step.

This completes the proof.
\end{proof}

\begin{proposition}
\label{PQPthExtEch}
Let $\boldsymbol{\mathcal{A}}^{(1)}$ be a first-order rewriting system and let $\mathrm{id}^{\boldsymbol{\mathcal{A}}^{(1)}}$ be its first-order identity morphism. Thus, $\mathrm{ech}^{(1,\mathcal{A}^{(1)})@}_{\boldsymbol{\mathcal{A}}^{(1)}} = \mathrm{id}^{[\mathrm{Pth}_{\boldsymbol{\mathcal{A}}^{(1)}}]}$.
\end{proposition}

\begin{proof}
Note that the following chain of equalities holds
\allowdisplaybreaks
\begin{align*}
\mathrm{id}^{[\mathbf{Pth}_{\boldsymbol{\mathcal{A}}^{(1)}}]}
\circ
\mathrm{pr}^{\mathrm{Ker}(\mathrm{CH}^{(1)})}_{\boldsymbol{\mathcal{A}}^{(1)}}
&=
\mathrm{ip}^{(1,X)@}_{\boldsymbol{\mathcal{A}}^{(1)}}
\circ
\mathrm{CH}^{[1]}_{\boldsymbol{\mathcal{A}}^{(1)}}
\circ
\mathrm{pr}^{\mathrm{Ker}(\mathrm{CH}^{(1)})}_{\boldsymbol{\mathcal{A}}^{(1)}}
\tag{1}
\\
&=
\mathrm{ip}^{(1,X)@}_{\boldsymbol{\mathcal{A}}^{(1)}}
\circ
\mathrm{pr}^{\Theta^{[1]}_{\boldsymbol{\mathcal{A}}^{(1)}}}
\circ
\mathrm{CH}^{(1)\mathrm{m}}_{\boldsymbol{\mathcal{A}}^{(1)}}
\circ
\mathrm{pr}^{\mathrm{Ker}(\mathrm{CH}^{(1)})}_{\boldsymbol{\mathcal{A}}^{(1)}}
\tag{2}
\\
&=
\mathrm{ip}^{(1,X)@}_{\boldsymbol{\mathcal{A}}^{(1)}}
\circ
\mathrm{pr}^{\Theta^{[1]}_{\boldsymbol{\mathcal{A}}^{(1)}}}
\circ
\mathrm{CH}^{(1)}_{\boldsymbol{\mathcal{A}}^{(1)}}
\tag{3}
\\
&=
\mathrm{pr}^{\mathrm{Ker}(\mathrm{CH}^{(1)})}_{\boldsymbol{\mathcal{A}}^{(1)}}
\circ
\mathrm{ip}^{(2,X)@}_{\boldsymbol{\mathcal{A}}^{(1)}}
\circ
\mathrm{CH}^{(1)}_{\boldsymbol{\mathcal{A}}^{(1)}}
\tag{4}
\end{align*}

The first equality follows from Theorem~\ref{TIso};
the second equality follows from the definition of the \(S\)-sorted mapping \(\mathrm{CH}^{[1]}_{\boldsymbol{\mathcal{A}}^{(1)}}\), introduced in Definition~\ref{DPTQCH};
the third equality follows from the definition of the monomorphic Curry-Howard mapping \(\mathrm{CH}^{(1)\mathrm{m}}_{\boldsymbol{\mathcal{A}}^{(1)}}\), introduced in Definition~\ref{DCHQuot};
finally, the last equality follows from the definition of the \(S\)-sorted mapping \(\mathrm{ip}^{(1,X)@}_{\boldsymbol{\mathcal{A}}^{(1)}}\), introduced in Definition~\ref{DPTQIp}.

Thus, by uniqueness of the first-order quotient path extension mapping introduced in Definition~\ref{DQPthExt}, \(\mathrm{ech}^{(1,\boldsymbol{\mathcal{A}}^{(1)})@}\) is equal to \(\mathrm{id}^{[\mathrm{Pth}_{\boldsymbol{\mathcal{A}}^{(1)}}]}\).
\end{proof}

Now we define the composition between first-order morphisms and show that it is not associative.

\begin{definition}
\label{DCompRws1}
Let $\mathbf{f}^{(1)}=(\varphi, c, (f^{(i)})_{i\in 2})$ be a first-order morphism from $\boldsymbol{\mathcal{A}}^{(1)}$ to $\boldsymbol{\mathcal{B}}^{(1)}$ and $\mathbf{g}^{(1)}=(\psi, d, (g^{(i)})_{i\in 2})$ a first-order morphism from $\boldsymbol{\mathcal{B}}^{(1)}$ to $\boldsymbol{\mathcal{C}}^{(1)}$. Its \emph{first-order composition}, from $\boldsymbol{\mathcal{A}}^{(1)}$ to $\boldsymbol{\mathcal{C}}^{(1)}$, is given by
$$
\mathbf{g}^{(1)} \circ \mathbf{f}^{(1)}
=
\left(
\mathbf{g}^{(0)} \circ \mathbf{f}^{(0)},
g^{(1)\flat}_{\varphi} \circ f^{(1)}
\right).
$$
Note that $\mathbf{g}^{(0)} \circ \mathbf{f}^{(0)}$ is the zeroth-order composition morphism, introduced in Definition~\ref{DCompRws0}, and $g^{(1)\flat}_{\varphi} \circ f^{(1)}$ is obtained from the following diagrams,
$$
\xymatrix{
\mathcal{B}^{(1)}
  \ar[r]^-{\mathrm{ech}^{(1,\mathcal{B}^{(1)})}_{\boldsymbol{\mathcal{B}}^{(1)}}}
  \ar[rd]_-{g^{(1)}}
&
\mathrm{Pth}_{\boldsymbol{\mathcal{B}}^{(1)}}
  \ar[d]^-{{g}^{(1)\flat}}
&&
\mathrm{Pth}_{\boldsymbol{\mathcal{B}}^{(1)}, \varphi}
  \ar[d]_-{g^{(1)\flat}_{\varphi}} 	
&
\mathcal{A^{(1)}}
  \ar[l]_-{f^{(1)}}
\\
& 
\mathrm{Pth}_{\boldsymbol{\mathcal{C}}^{(1)}, \psi}
&&
(\mathrm{Pth}_{\boldsymbol{\mathcal{C}}^{(1)}, \psi})_{\varphi}
\\
}
$$
being $g^{(1)\flat}$ the path extension mapping of $g^{(1)}$ introduced in Proposition~\ref{PDPthExt}.
\end{definition}


\begin{remark}
The composition of first-order morpshisms is not associative.
\end{remark}

We now define the notion of first-order equivalence of first-order morphisms. Thus, we will define the category $\mathsf{Rws}_{\mathfrak{d}}^{[1]}$ considering first-order rewriting system and first-order morphisms up to first-order equivalence, i.e., equivalence classes of first-order morphisms.

\begin{definition}
\label{DMorEqv}
Let $\mathbf{f}^{(1)}=(\varphi, c, (f^{(i)})_{i\in 2})$ and $\mathbf{g}^{(1)}=(\psi, d, (g^{(i)})_{i\in 2})$ be two first-order morphisms from $\boldsymbol{\mathcal{A}}^{(1)}$ to $\boldsymbol{\mathcal{B}}^{(1)}$. We will say that $\mathbf{f}^{(1)}$ and $\mathbf{g}^{(1)}$ are \emph{first-order equivalent}, written $\mathbf{f}^{(1)} \cong^{(1)} \mathbf{g}^{(1)}$, if 
\begin{enumerate}
\item
$(\varphi, c, f^{(0)}) = (\psi, d, g^{(0)})$, i.e., $\varphi = \psi$, $c=d$ and $f^{(0)}=g^{(0)}$; and
\item
$\mathrm{pr}^{\mathrm{Ker}(\mathrm{CH}^{(1)})}_{\boldsymbol{\mathcal{B}}^{(1)}, \varphi} \circ f^{(1)\flat} = \mathrm{pr}^{\mathrm{Ker}(\mathrm{CH}^{(1)})}_{\boldsymbol{\mathcal{B}}^{(1)}, \varphi} \circ g^{(1)\flat}$. That is, for every sort $s$ in $S$ and every path $\mathfrak{P}$ in $\mathrm{Pth}_{\boldsymbol{\mathcal{A}}, s}$,
$$
\left[
f_{s}^{(1)\flat}\left(
\mathfrak{P}
\right)
\right]_{\varphi(s)}
=
\left[
g_{s}^{(1)\flat}\left(
\mathfrak{P}
\right)
\right]_{\varphi(s)}.
$$
\end{enumerate}
Note that $\cong^{(1)}$ is an equivalence relation. Therefore, to simplify notation, we will denote by $[\mathbf{f}^{(1)}]$ the equivalence class $[\mathbf{f}^{(1)}]_{\cong^{(1)}}$.
\end{definition}

In the next proposition we state the relation between the path-extension mapping of a composition of first-order morphisms and the path-extension mappings of each of the composites. Let us say that, contrary to what happens with the homomorphic extension at level 0, it is not necessarily true that $(g^{(1)\flat}_{\varphi} \circ  f^{(1)} )^{\flat}$ is equal to $g^{(1)\flat}_{\varphi} \circ f^{(1)\flat}$.

\begin{proposition}
\label{PPthExtComp}
Let $\mathbf{f}^{(1)}=(\varphi, c, (f^{(i)})_{i\in 2})$ be a first-order morphism from $\boldsymbol{\mathcal{A}}^{(1)}$ to $\boldsymbol{\mathcal{B}}^{(1)}$ and $\mathbf{g}^{(1)}=(\psi, d, (g^{(i)})_{i\in 2})$ a first-order morphism from $\boldsymbol{\mathcal{B}}^{(1)}$ to $\boldsymbol{\mathcal{C}}^{(1)}$. Then
$$
\mathrm{pr}^{\mathrm{Ker}(\mathrm{CH}^{(1)})}_{\boldsymbol{\mathcal{C}}^{(1)}, \psi \circ \varphi}
\circ
\left(
g^{(1)\flat}_{\varphi} \circ 
f^{(1)}
\right)^{\flat}
=
\mathrm{pr}^{\mathrm{Ker}(\mathrm{CH}^{(1)})}_{\boldsymbol{\mathcal{C}}^{(1)}, \psi \circ \varphi}
\circ
g^{(1)\flat}_{\varphi} \circ
f^{(1)\flat}
$$
The reader is advised to consult the diagram of Figure~\ref{FPthExtComp}.

\begin{figure}
$$
\xymatrix{
\mathcal{A}^{(1)}
	\ar[r]^-{\mathrm{ech}^{(1,\boldsymbol{\mathcal{A}}^{(1)})}}
	\ar[rdd]_-{g^{(1)\flat}_{\varphi} \circ f^{(1)}}
	&
\mathrm{Pth}_{\boldsymbol{\mathcal{A}}^{(1)}}
	\ar[dd]^-{\left(g^{(1)\flat}_{\varphi} \circ f^{(1)}\right)^{\flat}}
	&&
\mathrm{Pth}_{\boldsymbol{\mathcal{A}}^{(1)}}
	\ar[d]_-{f^{(1)\flat}}
	&
\mathcal{A}^{(1)}
	\ar[l]_-{\mathrm{ech}^{(1,\boldsymbol{\mathcal{A}}^{(1)})}}
	\ar[ld]^-{f^{(1)}}
\\
	&&&
\mathrm{Pth}_{\boldsymbol{\mathcal{B}}^{(1)}, \varphi}
	\ar[d]_-{g^{(1)\flat}_{\varphi}}
\\
	&
\mathrm{Pth}_{\boldsymbol{\mathcal{C}}^{(1)}, \psi\circ\varphi}
	\ar[d]^-{\mathrm{pr}^{\mathrm{Ker}(\mathrm{CH}^{(1)})}_{\boldsymbol{\mathcal{C}}^{(1)}, \psi\circ\varphi}}
	&&
\mathrm{Pth}_{\boldsymbol{\mathcal{C}}^{(1)}, \psi\circ\varphi}
	\ar[d]_-{\mathrm{pr}^{\mathrm{Ker}(\mathrm{CH}^{(1)})}_{\boldsymbol{\mathcal{C}}^{(1)}, \psi\circ\varphi}}
\\
	&
[\mathrm{Pth}_{\boldsymbol{\mathcal{C}}^{(1)}}]_{\psi\circ\varphi}
	&&
[\mathrm{Pth}_{\boldsymbol{\mathcal{C}}^{(1)}}]_{\psi\circ\varphi}
}
$$
\caption{Relation between the path extension mapping of a composition and the corresponding path extension mappings.}
\label{FPthExtComp}
\end{figure}
\end{proposition}

\begin{proof}
Let $s$ in $S$ be a sort and $\mathfrak{P}$ a path in $\mathrm{Pth}_{\boldsymbol{\mathcal{A}}^{(1)}, s}$. We will prove that
$$
\left[
\left(g^{(1)\flat}_{\varphi} \circ f^{(1)}\right)^{\flat}_{s}\left(
\mathfrak{P}
\right)
\right]_{\psi(\varphi(s))}
=
\left[
g^{(1)\flat}_{\varphi(s)} \circ f^{(1)\flat}_{s}\left(
\mathfrak{P}
\right)
\right]_{\psi(\varphi(s))}
$$
by Artinian induction on $(\coprod \mathrm{Pth}_{\boldsymbol{\mathcal{A}}^{(1)}}, \leq_{\mathbf{Pth}_{\boldsymbol{\mathcal{A}}^{(1)}}})$.

{\sffamily Base step of the Artinian induction.}

Let $(\mathfrak{P}, s)$ be a minimal element of $(\coprod \mathrm{Pth}_{\boldsymbol{\mathcal{A}}^{(1)}}, \leq_{\mathbf{Pth}_{\boldsymbol{\mathcal{A}}^{(1)}}})$. Then, by Proposition~\ref{PMinimal}, the path $\mathfrak{P}$ is either~(1) an $(1,0)$-identity path or~(2) an echelon.

If~(1), i.e., if $\mathfrak{P}$ is a $(1,0)$-identity path, then $\mathfrak{P}=\mathrm{ip}_{\boldsymbol{\mathcal{A}}^{(1)},s}^{(1,0)\sharp}(P)$ for some term $P \in \mathrm{T}_{\Sigma}(X)$. Then the following chain of equalities holds
\allowdisplaybreaks
\begin{flalign*}
\left(
g^{(1)\flat}_{\varphi}
\circ
f^{(1)}
\right)^{\flat}_{s}\left(
\mathfrak{P}
\right)
&=
\left(
g^{(1)\flat}_{\varphi}
\circ
f^{(1)}
\right)^{\flat}_{s}\left(
\mathrm{ip}_{\boldsymbol{\mathcal{A}}^{(1)},s}^{(1,0])\sharp}\left(
P
\right)
\right)
\tag{1}
\\
&=
\mathrm{ip}_{\boldsymbol{\mathcal{C}}^{(1)}, \psi(\varphi(s))}^{(1,0)\sharp} \left(
\left(
g^{(0)\flat}_{\varphi} \circ f^{(0)}
\right)_{s}^{\sharp} \left(
P
\right)
\right)
\tag{2}
\\
&=
\mathrm{ip}_{\boldsymbol{\mathcal{C}}^{(1)}, \psi(\varphi(s))}^{(1,0)\sharp} \left(
g^{(0)\sharp}_{\varphi(s)} \left(
f_{s}^{(0)\sharp} \left(
P
\right) 
\right)
\right)
\tag{3}
\\
&=
g_{\varphi(s)}^{(1)\flat} \left(
\mathrm{ip}_{\boldsymbol{\mathcal{B}}^{(1)}, \varphi(s)}^{(1,0)\sharp} \left(
f_{s}^{(0)\sharp} \left(
\left[P\right]_{s}
\right) 
\right)
\right)
\tag{4}
\\
&=
g^{(1)\flat}_{\varphi(s)} \left(
f_{s}^{(1)\flat} \left(
\mathrm{ip}_{\boldsymbol{\mathcal{A}}^{(1)}, s}^{(1,0)\sharp} \left(
P
\right) 
\right)
\right)
\tag{5}
\\
&=
g^{(1)\flat}_{\varphi(s)} \left(
f_{s}^{(1)\flat} \left(
\mathfrak{P}
\right)
\right).
\tag{6}
\end{flalign*}

The first equality unravels the definition of the path $\mathfrak{P}$;
the second equality follows from Proposition~\ref{PPthExt};
the third equality follows from Remark~\ref{RCompExt} and unraveling the $s$-th component of the composition of $S$-sorted mappings;
the fourth and fifth equalities follow from Proposition~\ref{PPthExt};
finally, the last equality recovers the definition of the path $\mathfrak{P}$.

Therefore, their $\mathrm{Ker}(\mathrm{CH}^{(1)}_{\boldsymbol{\mathcal{C}}^{(1)}})$-classes coincide.

This completes case (1).

If~(2), i.e., if $\mathfrak{P}$ is an echelon associated to a rewrite rule $\mathfrak{p}=(M, N)$, that is, if $\mathfrak{P}$ has the form
$$
\xymatrix@C=55pt{
\mathfrak{P}: M
\ar@{->}[r]^-{\text{\Small{$(\mathfrak{p},\mathrm{id}^{\mathrm{T}_{\Sigma}(X)_{s}})$}}}
&
N
},
$$

Then the following chain of equalities holds
\begin{flalign*}
\left(
g^{(1)\flat}_{\varphi}
\circ
f^{(1)}
\right)^{\flat}_{s}\left(
\mathfrak{P}
\right)
&=
\left(
g^{(1)\flat}_{\varphi}
\circ
f^{(1)}
\right)^{\flat}_{s}\left(
\mathrm{ech}_{\boldsymbol{\mathcal{A}}^{(1)}, s}^{(1,\mathcal{A}^{(1)})}\left(
\mathfrak{p}
\right)
\right)
\tag{1}
\\
&=
\left(
g^{(1)\flat}_{\varphi}
\circ
f^{(1)}
\right)_{s}\left(
\mathfrak{p}
\right)
\tag{2}
\\
&=
g^{(1)\flat}_{\varphi(s)} \left(
f_{s}^{(1)} \left(
\mathfrak{p}
\right)
\right)
\tag{3}
\\
&=
g^{(1)\flat}_{\varphi(s)} \left(
f_{s}^{(1)\flat} \left(
\mathrm{ech}_{\boldsymbol{\mathcal{A}}^{(1)}, s}^{(1,\mathcal{A}^{(1)})}\left(
\mathfrak{p}
\right)
\right)
\right)
\tag{4}
\\
&=
g^{(1)\flat}_{\varphi(s)} \left(
f_{s}^{(1)\flat} \left(
\mathfrak{P}
\right)
\right).
\tag{5}
\end{flalign*}

The first equality unravels the definition of the path $\mathfrak{P}$;
the second equality follows from Proposition~\ref{PPthExt}; 
%
the third equality unravels the definition of the $s$-th component of the composition of $S$-sorted mappings;
the fourth equality follows from Proposition~\ref{PPthExt};
%
finally, the last equality recovers the definition of the path $\mathfrak{P}$.

Therefore, their $\mathrm{Ker}(\mathrm{CH}^{(1)}_{\boldsymbol{\mathcal{C}}^{(1)}})$-classes coincide.

This completes case (2).

This concludes the base step of the Artinian induction.

{\sffamily Inductive step of the Artinian induction.}

Let $(\mathfrak{P},s)$ be a non-minimal element of $(\coprod\mathrm{Pth}_{\boldsymbol{\mathcal{A}}^{(1)}}, \leq_{\mathbf{Pth}_{\boldsymbol{\mathcal{A}}^{(1)}}})$. We can assume that $\mathfrak{P}$ is a not a $(1,0)$-identity path, since for those paths the desired equality has already been proven. Let us suppose that, for every sort $t\in S$ and every path $\mathfrak{Q}\in\mathrm{Pth}_{\boldsymbol{\mathcal{A}}^{(1)},t}$, if $(\mathfrak{Q},t)<_{\mathbf{Pth}_{\boldsymbol{\mathcal{A}}^{(1)}}}(\mathfrak{P},s)$, then the equality
$$
\left[
\left(g^{(1)\flat}_{\varphi} \circ f^{(1)}\right)^{\flat}_{t}\left(
\mathfrak{Q}
\right)
\right]_{\psi(\varphi(t))}
=
\left[
g^{(1)\flat}_{\varphi(t)} \circ f^{(1)\flat}_{t}\left(
\mathfrak{Q}
\right)
\right]_{\psi(\varphi(t))}
$$
holds.

By Lemma~\ref{LOrdI}, we have that $\mathfrak{P}$ is either~(1) a path of length strictly greater than one containing at least one first-order echelon or~(2) an
echelonless path.

If~(1), i.e., if $\mathfrak{P}$ is a path of length strictly greater than one containing at least one first-order echelon, then let $i\in \bb{\mathfrak{P}}$ be the first index for which the one-step subpath $\mathfrak{P}^{i,i}$ of $\mathfrak{P}$ is an echelon. We consider different cases for $i$ according to the cases presented in Definition~\ref{DOrd}.

If~$i=0$, we have that the pairs $(\mathfrak{P}^{0,0}, s)$ and $(\mathfrak{P}^{1,\bb{\mathfrak{P}}-1}, s)$ $\prec_{\mathbf{Pth}_{\boldsymbol{\mathcal{A}}^{(1)}}}$-precede the pair $(\mathfrak{P}, s)$. Therefore, by induction, the equalities
$$
\left[
\left(g^{(1)\flat}_{\varphi} \circ f^{(1)}\right)^{\flat}_{s}\left(
\mathfrak{P}^{0,0}
\right)
\right]_{\psi(\varphi(s))}
=
\left[
g^{(1)\flat}_{\varphi(s)} \circ f^{(1)\flat}_{s}\left(
\mathfrak{P}^{0,0}
\right)
\right]_{\psi(\varphi(s))}
$$
$$
\left[
\left(g^{(1)\flat}_{\varphi} \circ f^{(1)}\right)^{\flat}_{s}\left(
\mathfrak{P}^{1,\bb{\mathfrak{P}}-1}
\right)
\right]_{\psi(\varphi(s))}
=
\left[
g^{(1)\flat}_{\varphi(s)} \circ f^{(1)\flat}_{s}\left(
\mathfrak{P}^{1,\bb{\mathfrak{P}}-1}
\right)
\right]_{\psi(\varphi(s))}
$$
hold. Then, the following chain of equalities holds
\begin{flushleft}
$
\left[
\left(
g^{(1)\flat}_{\varphi} \circ f^{(1)}
\right)_{s}^{\flat} \left(
\mathfrak{P}
\right)
\right]_{\psi(\varphi(s))}
$
\allowdisplaybreaks
\begin{align*}
&=
\left[
\left(
g^{(1)\flat}_{\varphi} \circ f^{(1)}
\right)_{s}^{\flat} \left(
\mathfrak{P}^{1, \bb{\mathfrak{P}}-1} \circ_{s}^{0\mathbf{Pth}_{\boldsymbol{\mathcal{A}}^{(1)}}} \mathfrak{P}^{(1)0,0}
\right)
\right]_{\psi(\varphi(s))}
\tag{1}
\\
&=
\left[
\left(
g^{(1)\flat}_{\varphi} \circ f^{(1)}
\right)_{s}^{\flat} \left(
\mathfrak{P}^{1, \bb{\mathfrak{P}}-1}
\right)
\right]_{\psi(\varphi(s))}
\\
&\hspace{3.5cm}
\circ_{s}^{0\left[\mathbf{Pth}_{\boldsymbol{\mathcal{C}}^{(1)}}^{\mathbf{g}^{(1)}\circ\mathbf{f}^{(1)}}\right]}
\left[
\left(
g^{(1)\flat}_{\varphi} \circ f^{(1)}
\right)_{s}^{\flat} \left(
\mathfrak{P}^{0,0}
\right)
\right]_{\psi(\varphi(s))}
\tag{2}
\\
&=
\left[
\left(g^{(1)\flat}_{\varphi} \circ f^{(1)\flat}\right)_{s} \left(
\mathfrak{P}^{1, \bb{\mathfrak{P}}-1}
\right)
\right]_{\psi(\varphi(s))}
\\
&\hspace{3.5cm}
\circ_{s}^{0\left[\mathbf{Pth}_{\boldsymbol{\mathcal{C}}^{(1)}}^{\mathbf{g}^{(1)}\circ\mathbf{f}^{(1)}}\right]}
\left[
\left(g^{(1)\flat}_{\varphi} \circ f^{(1)\flat}\right)_{s} \left(
\mathfrak{P}^{0,0}
\right)
\right]_{\psi(\varphi(s))}
\tag{3}
\\
&=
\left[
g_{\varphi(s)}^{(1)\flat} \left(
f_{s}^{(1)\flat} \left(
\mathfrak{P}^{1, \bb{\mathfrak{P}}-1}
\right)
\right)
\right]_{\psi(\varphi(s))}
\\
&\hspace{3.5cm}
\circ_{\varphi(s)}^{0[\mathbf{Pth}_{\boldsymbol{\mathcal{C}}^{(1)}}^{\mathbf{g}^{(1)}}]}
\left[
g_{\varphi(s)}^{(1)\flat} \left(
f_{s}^{(1)\flat} \left(
\mathfrak{P}^{0,0}
\right)
\right)
\right]_{\psi(\varphi(s))}
\tag{4}
\\
&=
\left[
g_{\varphi(s)}^{(1)\flat} \left(
f_{s}^{(1)\flat} \left(
\mathfrak{P}^{1, \bb{\mathfrak{P}}-1}
\right)
\circ_{\varphi(s)}^{0\mathbf{Pth}_{\boldsymbol{\mathcal{B}}^{(1)}}}
f_{s}^{(1)\flat} \left(
\mathfrak{P}^{0,0}
\right)
\right)
\right]_{\psi(\varphi(s))}
\tag{5}
\\
&=
\left[
g_{\varphi(s)}^{(1)\flat} \left(
f_{s}^{(1)\flat} \left(
\mathfrak{P}^{1, \bb{\mathfrak{P}}-1}
\circ_{s}^{0\mathbf{Pth}_{\boldsymbol{\mathcal{A}}^{(1)}}}
\mathfrak{P}^{0,0}
\right)
\right)
\right]_{\psi(\varphi(s))}
\tag{6}
\\
&=
\left[
g_{\varphi(s)}^{(1)\flat} \left(
f_{s}^{(1)\flat} \left(
\mathfrak{P}
\right)
\right)
\right]_{\psi(\varphi(s))}.
\tag{7}
\end{align*}
\end{flushleft}

The first equality follows from the fact that $\mathfrak{P} = \mathfrak{P}^{1, \bb{\mathfrak{P}}-1} \circ_{s}^{0\mathbf{Pth}_{\boldsymbol{\mathcal{A}}^{(1)}}} \mathfrak{P}^{0,0}$;
the second equality follows from Proposition~\ref{PHomPthExtKer};
the third equality follows by Artinian induction;
the fourth equality unravels the definition of the $s$-th component of the composition of $S$-sorted mappings and note that, for every sort $s$ in $S$,
$$
\circ_{s}^{0[\mathbf{Pth}_{\boldsymbol{\mathcal{C}}^{(1)}}^{\mathbf{g}^{(1)}\circ\mathbf{f}^{(1)}}]}
=
\circ_{\psi(\varphi(s))}^{0[\mathbf{Pth}_{\boldsymbol{\mathcal{C}}^{(1)}}]}
=
\circ_{\varphi(s)}^{0[\mathbf{Pth}_{\boldsymbol{\mathcal{C}}^{(1)}}^{\mathbf{g}^{(1)}}]};
$$
the fifth equality follows from Proposition~\ref{PHomPthExtKer};
the sixth equality recovers the definition of $f^{(1)\flat}$;
finally, the last equality recovers the definition of $\mathfrak{P}$;

If~$i>0$, we have that the pairs $(\mathfrak{P}^{0,i-1}, s)$ and $(\mathfrak{P}^{i,\bb{\mathfrak{P}}-1}, s)$ $\prec_{\mathbf{Pth}_{\boldsymbol{\mathcal{A}}^{(1)}}}$-precede the pair $(\mathfrak{P}, s)$. Therefore, by induction, the equalities
$$
\left[
\left(g^{(1)\flat}_{\varphi} \circ f^{(1)}\right)^{\flat}_{s}\left(
\mathfrak{P}^{0,i-1}
\right)
\right]_{\psi(\varphi(s))}
=
\left[
g^{(1)\flat}_{\varphi(s)} \circ f^{(1)\flat}_{s}\left(
\mathfrak{P}^{0,i-1}
\right)
\right]_{\psi(\varphi(s))}
$$
$$
\left[
\left(g^{(1)\flat}_{\varphi} \circ f^{(1)}\right)^{\flat}_{s}\left(
\mathfrak{P}^{i,\bb{\mathfrak{P}}-1}
\right)
\right]_{\psi(\varphi(s))}
=
\left[
g^{(1)\flat}_{\varphi(s)} \circ f^{(1)\flat}_{s}\left(
\mathfrak{P}^{i,\bb{\mathfrak{P}}-1}
\right)
\right]_{\psi(\varphi(s))}
$$
hold. Then, the following chain of equalities holds
\begin{flushleft}
$
\left[
\left(
g^{(1)\flat}_{\varphi} \circ f^{(1)}
\right)_{s}^{\flat} \left(
\mathfrak{P}
\right)
\right]_{\psi(\varphi(s))}
$
\allowdisplaybreaks
\begin{align*}
&=
\left[
\left(
g^{(1)\flat}_{\varphi} \circ f^{(1)}
\right)_{s}^{\flat} \left(
\mathfrak{P}^{i, \bb{\mathfrak{P}}-1} \circ_{s}^{0\mathbf{Pth}_{\boldsymbol{\mathcal{A}}^{(1)}}}
\mathfrak{P}^{0,i-1}
\right)
\right]_{\psi(\varphi(s))}
\tag{1}
\\
&=
\left[
\left(
g^{(1)\flat}_{\varphi} \circ f^{(1)}
\right)_{s}^{\flat} \left(
\mathfrak{P}^{i, \bb{\mathfrak{P}}-1}
\right)
\right]_{\psi(\varphi(s))}
\\
&\hspace{3.5cm}
\circ_{s}^{0\left[\mathbf{Pth}_{\boldsymbol{\mathcal{C}}^{(1)}}^{\mathbf{g}^{(1)}\circ\mathbf{f}^{(1)}}\right]}
\left[
\left(
g^{(1)\flat}_{\varphi} \circ f^{(1)}
\right)_{s}^{\flat} \left(
\mathfrak{P}^{0,i-1}
\right)
\right]_{\psi(\varphi(s))}
\tag{2}
\\
&=
\left[
\left(g^{(1)\flat}_{\varphi} \circ f^{(1)\flat}\right)_{s} \left(
\mathfrak{P}^{i, \bb{\mathfrak{P}}-1}
\right)
\right]_{\psi(\varphi(s))}
\\
&\hspace{3.5cm}
\circ_{s}^{0\left[\mathbf{Pth}_{\boldsymbol{\mathcal{C}}^{(1)}}^{\mathbf{g}^{(1)}\circ\mathbf{f}^{(1)}}\right]}
\left[
\left(g^{(1)\flat}_{\varphi} \circ f^{(1)\flat}\right)_{s} \left(
\mathfrak{P}^{0,i-1}
\right)
\right]_{\psi(\varphi(s))}
\tag{3}
\\
&=
\left[
g_{\varphi(s)}^{(1)\flat} \left(
f_{s}^{(1)\flat} \left(
\mathfrak{P}^{i, \bb{\mathfrak{P}}-1}
\right)
\right)
\right]_{\psi(\varphi(s))}
\\
&\hspace{3.5cm}
\circ_{\varphi(s)}^{0[\mathbf{Pth}_{\boldsymbol{\mathcal{C}}^{(1)}}^{\mathbf{g}^{(1)}}]}
\left[
g_{\varphi(s)}^{(1)\flat} \left(
f_{s}^{(1)\flat} \left(
\mathfrak{P}^{0,i-1}
\right)
\right)
\right]_{\psi(\varphi(s))}
\tag{4}
\\
&=
\left[
g_{\varphi(s)}^{(1)\flat} \left(
f_{s}^{(1)\flat} \left(
\mathfrak{P}^{i, \bb{\mathfrak{P}}-1}
\right)
\circ_{\varphi(s)}^{0\mathbf{Pth}_{\boldsymbol{\mathcal{B}}^{(1)}}}
f_{s}^{(1)\flat} \left(
\mathfrak{P}^{0,i-1}
\right)
\right)
\right]_{\psi(\varphi(s))}
\tag{5}
\\
&=
\left[
g_{\varphi(s)}^{(1)\flat} \left(
f_{s}^{(1)\flat} \left(
\mathfrak{P}^{i, \bb{\mathfrak{P}}-1}
\circ_{s}^{0\mathbf{Pth}_{\boldsymbol{\mathcal{A}}^{(1)}}}
\mathfrak{P}^{0,i-1}
\right)
\right)
\right]_{\psi(\varphi(s))}
\tag{6}
\\
&=
\left[
g_{\varphi(s)}^{(1)\flat} \left(
f_{s}^{(1)\flat} \left(
\mathfrak{P}
\right)
\right)
\right]_{\psi(\varphi(s))}.
\tag{7}
\end{align*}
\end{flushleft}

The first equality follows from the fact that $\mathfrak{P} = \mathfrak{P}^{i, \bb{\mathfrak{P}}-1} \circ_{s}^{0\mathbf{Pth}_{\boldsymbol{\mathcal{A}}^{(1)}}} \mathfrak{P}^{0,i-1}$;
the second equality follows from Proposition~\ref{PHomPthExtKer};
the third equality follows by Artinian induction;
the fourth equality unravels the definition of the $s$-th component of the composition of $S$-sorted mappings and note that, for every sort $s$ in $S$,
$$
\circ_{s}^{0[\mathbf{Pth}_{\boldsymbol{\mathcal{C}}^{(1)}}^{\mathbf{g}^{(1)}\circ\mathbf{f}^{(1)}}]}
=
\circ_{\psi(\varphi(s))}^{0[\mathbf{Pth}_{\boldsymbol{\mathcal{C}}^{(1)}}]}
=
\circ_{\varphi(s)}^{0[\mathbf{Pth}_{\boldsymbol{\mathcal{C}}^{(1)}}^{\mathbf{g}^{(1)}}]};
$$
the fifth equality follows from Proposition~\ref{PHomPthExtKer};
the sixth equality recovers the definition of $f^{(1)\flat}$;
finally, the last equality recovers the definition of $\mathfrak{P}$;

Case (1) follows.

If~(2), i.e., if $\mathfrak{P}$ is an echelonless path in $\mathrm{Pth}_{\boldsymbol{\mathcal{A}}^{(1)},s}$, then the conditions for the first-order extraction algorithm, that is, Lemma~\ref{LPthExtract}, are fulfilled. Then there exists a unique word $\mathbf{s} \in S^{\star} - \{\lambda\}$ and a unique operation symbol $\sigma \in \Sigma_{\mathbf{s}, s}$ associated to $\mathfrak{P}$. Let $(\mathfrak{P}_{j})_{j \in \bb{\mathbf{s}}}$ be the family of paths in $\mathrm{Pth}_{\boldsymbol{\mathcal{A}}^{(1)}, \mathbf{s}}$, which, in virtue of Lemma~\ref{LPthExtract}, we can extract from $\mathfrak{P}$. Note that, for every $j \in \bb{\mathbf{s}}$, the pair $(\mathfrak{P}_{j}, s_{j})$ $\prec_{\mathbf{Pth}_{\boldsymbol{\mathcal{A}}^{(1)}}}$- precede the pair $(\mathfrak{P}, s)$. Therefore, by induction, for every $j \in \bb{\mathbf{s}}$, the equality
$$
\left[
\left(g^{(1)\flat}_{\varphi} \circ f^{(1)}\right)^{\flat}_{s_{j}}\left(
\mathfrak{P}_{j}
\right)
\right]_{\psi(\varphi(s_{j}))}
=
\left[
g^{(1)\flat}_{\varphi(s_{j})} \circ f^{(1)\flat}_{s_{j}}\left(
\mathfrak{P}_{j}
\right)
\right]_{\psi(\varphi(s_{j}))}
$$
hold.

The following chain of equalities holds
\begin{flushleft}
$
\left[
\left(
g^{(1)\flat}_{\varphi} \circ f^{(1)}
\right)_{s}^{\flat} \left(
\mathfrak{P}
\right)
\right]_{\psi(\varphi(s))}
$
\allowdisplaybreaks
\begin{align*}
&=
\left[
\sigma^{\mathbf{Pth}_{\boldsymbol{\mathcal{C}}^{(1)}}^{\mathbf{g}^{(1)} \circ \mathbf{f}^{(1)}}} \left(
\left(
\left(
g_{\varphi}^{(1)\flat} \circ
f^{(1)}
\right)_{s_{j}}^{\flat}
\left(
\mathfrak{P}_{j}
\right)
\right)_{j \in \bb{\mathbf{s}}}
\right)
\right]_{\psi(\varphi(s))}
\tag{1}
\\
&=
\sigma^{[\mathbf{Pth}_{\boldsymbol{\mathcal{C}}^{(1)}}^{\mathbf{g}^{(1)} \circ \mathbf{f}^{(1)}}]} \left(
\left(
\left[
\left(
g_{\varphi}^{(1)\flat} \circ
f^{(1)}
\right)_{s_{j}}^{\flat}
\left(
\mathfrak{P}_{j}
\right)
\right]_{\psi(\varphi(s_{j}))}
\right)_{j \in \bb{\mathbf{s}}}
\right)
\tag{2}
\\
&=
\sigma^{[\mathbf{Pth}_{\boldsymbol{\mathcal{C}}^{(1)}}^{\mathbf{g}^{(1)} \circ \mathbf{f}^{(1)}}]} \left(
\left(
\left[
g^{(1)\flat}_{\varphi(s_{j})} \circ 
f^{(1)\flat}_{s_{j}}\left(
\mathfrak{P}_{j}
\right)
\right]_{\psi(\varphi(s_{j}))}
\right)_{j \in \bb{\mathbf{s}}}
\right)
\tag{3}
\\
&=
\left[
\sigma^{\mathbf{Pth}_{\boldsymbol{\mathcal{C}}^{(1)}}^{\mathbf{g}^{(1)} \circ \mathbf{f}^{(1)}}} \left(
\left(
g_{\varphi(s_{j)}}^{(1)\flat} \circ
f^{(1)\flat}_{s_{j}}
\left(
\mathfrak{P}_{j}
\right)
\right)_{j \in \bb{\mathbf{s}}}
\right)
\right]_{\psi(\varphi(s))}
\tag{4}
\\
&=
\left[
\sigma^{\mathbf{d}_{\mathfrak{d}}^{\ast}(\mathbf{c}_{\mathfrak{d}}^{\ast}(\mathbf{Pth}_{\boldsymbol{\mathcal{C}}^{(1)}}))} \left(
\left(
g_{\varphi(s_{j})}^{(1)\flat} \circ
f^{(1)\flat}_{s_{j}}
\left(
\mathfrak{P}_{j}
\right)
\right)_{j \in \bb{\mathbf{s}}}
\right)
\right]_{\psi(\varphi(s))}
\tag{5}
\\
&=
\left[
g_{\varphi(s)}^{(1)\flat} \left(
\sigma^{\mathbf{c}_{\mathfrak{d}}^{\ast}(\mathbf{Pth}_{\boldsymbol{\mathcal{B}}^{(1)}})} \left(
\left(
f^{(1)\flat}_{s_{j}}
\left(
\mathfrak{P}_{j}
\right)
\right)_{j \in \bb{\mathbf{s}}}
\right)
\right)
\right]_{\psi(\varphi(s))}
\tag{6}
\\
&=
\left[
g^{(1)\flat}_{\varphi(s)} \left( 
f^{(1)\flat}_{s}\left(
\mathfrak{P}
\right)
\right)
\right]_{\psi(\varphi(s))}
\tag{7}
\end{align*}
\end{flushleft}

The first equality unravels the definition of the path extension mapping introduced in Proposition~\ref{PPthExt};
%
the second equality recovers the definition of the interpretation of $\sigma$ in the partial $\Sigma^{\boldsymbol{\mathcal{A}}^{(1)}}$-algebra $[ \mathbf{Pth}_{\boldsymbol{\mathcal{C}}^{(1)}}^{\mathbf{g}^{(1)} \circ \mathbf{f}^{(1)}} ]$ introduced in Proposition~\ref{PQPthBCatAlg};
the third equality follows by Artinian induction;
the fourth equality unravels the definition of the interpretation of $\sigma$ in the partial $\Sigma^{\boldsymbol{\mathcal{A}}^{(1)}}$-algebra $[ \mathbf{Pth}_{\boldsymbol{\mathcal{C}}^{(1)}}^{\mathbf{g}^{(1)} \circ \mathbf{f}^{(1)}} ]$ introduced in Proposition~\ref{PQPthBCatAlg};
the fifth equality unravels the definition of the interpretation of $\sigma$ in the partial $\Sigma^{\boldsymbol{\mathcal{A}}^{(1)}}$-algebra $\mathbf{Pth}_{\boldsymbol{\mathcal{C}}^{(1)}}^{\mathbf{g}^{(1)} \circ \mathbf{f}^{(1)}}$ introduced in Proposition~\ref{PPthBCatAlg};
the sixth equality follows from the fact that, according to Proposition~\ref{PPthExtHom}, $g^{(1)\flat}$ is a $\Lambda$-homomorphism, thus, according to Proposition~\ref{PFunSig}, $\mathbf{c}^{\ast}_{\mathfrak{d}}(g^{(1)\flat})$ is a $\Sigma$-homomorphism;
finally, the last equality follows from the fact that, according to Proposition~\ref{PPthExtHom}, $f^{(1)\flat}$ is a $\Sigma$-homomorphism.

This completes case (2).

This completes the proof.
\end{proof}

In the next corollary, we give equivalent expressions for the homomorphic path extension mapping of a first-order composition of first-order morphisms in terms of the respective homomorphic path extension mappings of the composites.

\begin{corollary}\label{CHomPthExtComp}
Let $\mathbf{f}^{(1)}=(c, \varphi, (f^{(i)})_{i\in 2})$ be a first-order morphism from $\boldsymbol{\mathcal{A}}^{(1)}$ to $\boldsymbol{\mathcal{B}}^{(1)}$ and $\mathbf{g}^{(1)}=(d, \psi, (g^{(i)})_{i\in 2})$ a first-order morphism from $\boldsymbol{\mathcal{B}}^{(1)}$ to $\boldsymbol{\mathcal{C}}^{(1)}$. Then the following three mappings are equal:
\begin{enumerate}
\item
$\mathrm{pr}^{\mathrm{Ker}(\mathrm{CH}^{(1)})}_{\boldsymbol{\mathcal{C}}^{(1)}, \psi\circ\varphi} \circ \left(g^{(1)\flat}_{\varphi} \circ f^{(1)}\right)^{\flat}$;
\item
$\left(\mathrm{pr}^{\mathrm{Ker}(\mathrm{CH}^{(1)})}_{\boldsymbol{\mathcal{C}}^{(1)}, \psi} \circ g^{(1)\flat}\right)_{\varphi} \circ f^{(1)\flat}$;
\item
$g^{(1)@}_{\varphi} \circ \left(\mathrm{pr}^{\mathrm{Ker}(\mathrm{CH}^{(1)})}_{\boldsymbol{\mathcal{B}}^{(1)}, \varphi} \circ f^{(1)\flat}\right)$.
\end{enumerate}
\end{corollary}

Finally, in the next corollary we give an expression for the quotient path extension mapping of the first-order composition of first-order morphisms in terms of the respective quotient path extension mappings of the composites.

\begin{corollary}\label{CQPthExtComp}
Let $\mathbf{f}^{(1)}=(\varphi, c, (f^{(i)})_{i\in 2})$ be a first-order morphism from $\boldsymbol{\mathcal{A}}^{(1)}$ to $\boldsymbol{\mathcal{B}}^{(1)}$ and $\mathbf{g}^{(1)}=(\psi, d, (g^{(i)})_{i\in 2})$ a first-order morphism from $\boldsymbol{\mathcal{B}}^{(1)}$ to $\boldsymbol{\mathcal{C}}^{(1)}$. Then 
$$
\left(
g^{(1)\flat}_{\varphi} \circ f^{(1)}
\right)^{@}
=
g^{[1]@}_{\varphi} \circ f^{[1]@}.
$$
\end{corollary}

We are now in a position where we can prove that the first-order rewriting systems and equivalence classes of first-order morphisms constitute a category. Before doing so, we show that the composition of classes of first-order morphisms is a well-defined binary operation and it does not depend on the choice of the representatives.

\begin{proposition}
\label{PCompRws1}
Let $\mathbf{f}^{(1)}$ and $\mathbf{f}'^{(1)}$ be first-order morphisms from $\boldsymbol{\mathcal{A}}^{(1)}$ to $\boldsymbol{\mathcal{B}}^{(1)}$ and $\mathbf{g}^{(1)}$ and $\mathbf{g}'^{(1)}$ first-order morphisms from $\boldsymbol{\mathcal{B}}^{(1)}$ to $\boldsymbol{\mathcal{C}}^{(1)}$. If $\mathbf{f}^{(1)} \cong^{(1)} \mathbf{f}'^{(1)}$ and $\mathbf{g}^{(1)} \cong^{(1)} \mathbf{g}'^{(1)}$, then
$
\mathbf{g}^{(1)} \circ \mathbf{f}^{(1)}
\cong^{(1)}
\mathbf{g}'^{(1)} \circ \mathbf{f}'^{(1)}.
$
Thus, the composition of their morphism classes $[\mathbf{g}^{(1)}] \circ [\mathbf{f}^{(1)}] = [\mathbf{g}^{(1)} \circ \mathbf{f}^{(1)}]$ is well-defined and does not depend on the representatives of the first-order morphism classes.
\end{proposition}

\begin{proof}
Let $\mathbf{f}^{(1)} = (\mathbf{f}^{(0)}, f^{(1)})$ and $\mathbf{f}'^{(1)} = (\mathbf{f}'^{(0)}, f'^{(1)})$. Since $\mathbf{f}^{(1)}$ and $\mathbf{f}'^{(1)}$ are first-order equivalent, following Definition~\ref{DMorEqv}, $\mathbf{f}^{(0)} = \mathbf{f}'^{(0)}$ and 
\begin{equation}
\mathrm{pr}^{\mathrm{Ker}(\mathrm{CH}^{(1)})}_{\boldsymbol{\mathcal{B}}^{(1)}, \varphi} \circ
f^{(1)\flat}
=
\mathrm{pr}^{\mathrm{Ker}(\mathrm{CH}^{(1)})}_{\boldsymbol{\mathcal{B}}^{(1)}, \varphi} \circ
f'^{(1)\flat}.
\tag{F}
\end{equation}
Moreover, let $\mathbf{g}^{(1)} = (\mathbf{g}^{(0)}, g^{(1)})$ and $\mathbf{g}'^{(1)} = (\mathbf{g}'^{(0)}, g'^{(1)})$. Since $\mathbf{g}^{(1)}$ and $\mathbf{g}'^{(1)}$ are first-order equivalent, following Definition~\ref{DMorEqv}, $\mathbf{g}^{(0)} = \mathbf{g}'^{(0)}$ and 
\begin{equation}
\mathrm{pr}^{\mathrm{Ker}(\mathrm{CH}^{(1)})}_{\boldsymbol{\mathcal{C}}^{(1)}, \psi}
\circ
g^{(1)\flat}
=
\mathrm{pr}^{\mathrm{Ker}(\mathrm{CH}^{(1)})}_{\boldsymbol{\mathcal{C}}^{(1)}, \psi}
\circ
g'^{(1)\flat}.
\tag{G}
\end{equation}

Following Definition~\ref{DCompRws1}, the respective compositions are given by
$$
\mathbf{g}^{(1)} \circ \mathbf{f}^{(1)}
=
\left(
\mathbf{g}^{(0)} \circ \mathbf{f}^{(0)},
g^{(1)\flat}_{\varphi} \circ f^{(1)}
\right)
$$
and
$$
\mathbf{g}'^{(0)} \circ \mathbf{f}'^{(0)}
=
\left(
\mathbf{g}'^{(1)} \circ \mathbf{f}'^{(1)},
g'^{(1)\flat}_{\varphi} \circ f'^{(1)}
\right).
$$

Note that, $\mathbf{g}^{(0)} \circ \mathbf{f}^{(0)}$ is equal to  $\mathbf{g}'^{(0)} \circ \mathbf{f}'^{(0)}$. Therefore, by Definition~\ref{DMorEqv}, to show that $\mathbf{g}^{(1)} \circ \mathbf{f}^{(1)}$ and $\mathbf{g}'^{(1)} \circ \mathbf{f}'^{(1)}$ are first-order equivalent, it suffices to show that
$$
\mathrm{pr}^{\mathrm{Ker}(\mathrm{CH}^{(1)})}_{\boldsymbol{\mathcal{C}}^{(1)}, \psi \circ \varphi}
\circ
\left(
g^{(1)\flat}_{\varphi} \circ 
f^{(1)}
\right)^{\flat}
=
\mathrm{pr}^{\mathrm{Ker}(\mathrm{CH}^{(1)})}_{\boldsymbol{\mathcal{C}}^{(1)}, \psi \circ \varphi}
\circ
\left(
g'^{(1)\flat}_{\varphi} \circ 
f'^{(1)}
\right)^{\flat}
$$

The following chain of equalities holds
\allowdisplaybreaks
\begin{align*}
\mathrm{pr}^{\mathrm{Ker}(\mathrm{CH}^{(1)})}_{\boldsymbol{\mathcal{C}}^{(1)}, \psi \circ \varphi}
\circ
\left(
g^{(1)\flat}_{\varphi} \circ 
f^{(1)}
\right)^{\flat}
&=
\mathrm{pr}^{\mathrm{Ker}(\mathrm{CH}^{(1)})}_{\boldsymbol{\mathcal{C}}^{(1)}, \psi \circ \varphi}
\circ
g^{(1)\flat}_{\varphi} \circ
f^{(1)\flat}
\tag{1}
\\
&=
\left(
\mathrm{pr}^{\mathrm{Ker}(\mathrm{CH}^{(1)})}_{\boldsymbol{\mathcal{C}}^{(1)}, \psi}
\right)_{\varphi} \circ
g^{(1)\flat}_{\varphi} \circ
f^{(1)\flat}
\tag{2}
\\
&=
\left(
\mathrm{pr}^{\mathrm{Ker}(\mathrm{CH}^{(1)})}_{\boldsymbol{\mathcal{C}}^{(1)}, \psi}
\circ
g^{(1)\flat}\right)_{\varphi} \circ
f^{(1)\flat}
\tag{3}
\\
&=
\left(
\mathrm{pr}^{\mathrm{Ker}(\mathrm{CH}^{(1)})}_{\boldsymbol{\mathcal{C}}^{(1)}, \psi}
\circ
g'^{(1)\flat}\right)_{\varphi} \circ
f^{(1)\flat}
\tag{4}
\\
&=
\left(
g'^{(1)@}
\circ
\mathrm{pr}^{\mathrm{Ker}(\mathrm{CH}^{(1)})}_{\boldsymbol{\mathcal{B}}^{(1)}}\right)_{\varphi} \circ
f^{(1)\flat}
\tag{5}
\\
&=
\left(
g'^{(1)@}
\right)_{\varphi}
\circ
\mathrm{pr}^{\mathrm{Ker}(\mathrm{CH}^{(1)})}_{\boldsymbol{\mathcal{B}}^{(1)}, \varphi} \circ
f^{(1)\flat}
\tag{6}
\\
&=
\left(
g'^{(1)@}
\right)_{\varphi}
\circ
\mathrm{pr}^{\mathrm{Ker}(\mathrm{CH}^{(1)})}_{\boldsymbol{\mathcal{B}}^{(1)}, \varphi} \circ
f'^{(1)\flat}
\tag{7}
\\
&=
\left(
g'^{(1)@}
\circ
\mathrm{pr}^{\mathrm{Ker}(\mathrm{CH}^{(1)})}_{\boldsymbol{\mathcal{B}}^{(1)}} \right)_{\varphi} \circ
f'^{(1)\flat}
\tag{8}
\\
&=
\left(
\mathrm{pr}^{\mathrm{Ker}(\mathrm{CH}^{(1)})}_{\boldsymbol{\mathcal{C}}^{(1)}, \psi}
\circ
g'^{(1)\flat}  \right)_{\varphi}
\circ
f'^{(1)\flat}
\tag{9}
\\
&=
\left(
\mathrm{pr}^{\mathrm{Ker}(\mathrm{CH}^{(1)})}_{\boldsymbol{\mathcal{C}}^{(1)}, \psi}
\right)_{\varphi}
\circ
g'^{(1)\flat}_{\varphi}
\circ
f'^{(1)\flat}
\tag{10}
\\
&=
\mathrm{pr}^{\mathrm{Ker}(\mathrm{CH}^{(1)})}_{\boldsymbol{\mathcal{C}}^{(1)}, \psi \circ \phi}
\circ
g'^{(1)\flat}_{\varphi}
\circ
f'^{(1)\flat}
\tag{11}
\\
&=
\mathrm{pr}^{\mathrm{Ker}(\mathrm{CH}^{(1)})}_{\boldsymbol{\mathcal{C}}^{(1)}, \psi \circ \phi}
\circ
\left(
g'^{(1)\flat}_{\varphi}
\circ
f'^{(1)}
\right)^{\flat}.
\tag{12}
\end{align*}

The first equality follows from Proposition~\ref{PPthExtComp};
the second equality follows from the fact that, according to Proposition~\ref{PMSetFunc}, $\mathrm{MSet}$ is a contravariant functor;
the third equality follows from the fact that, according to Proposition~\ref{PDeltaPhiFunc}, $\Delta_{\varphi}$ is a covariant functor;
the fourth equality follows from Equation~(G);
the fifth equality follows from the definition of $g'^{(1)@}$, introduced in Definition~\ref{DQPthExt};
the sixth equality follows from the fact that, according to Proposition~\ref{PDeltaPhiFunc}, $\Delta_{\varphi}$ is a covariant functor;
the seventh equality follows from Equation~(F);
the eighth equality follows from the fact that, according to Proposition~\ref{PDeltaPhiFunc}, $\Delta_{\varphi}$ is a covariant functor;
the ninth equality follows from the definition of $g'^{(1)@}$, introduced in Definition~\ref{DQPthExt};
the tenth equality follows from the fact that, according to Proposition~\ref{PDeltaPhiFunc}, $\Delta_{\varphi}$ is a covariant functor;
the eleventh equality follows from the fact that, according to Proposition~\ref{PMSetFunc}, $\mathrm{MSet}$ is a contravariant functor;
finally, the last equality follows from Proposition~\ref{PPthExtComp}.

This completes the proof.
\end{proof}

We now show that first-order many-sorted rewriting systems and classes of first-order morphisms between them constitute a category.

\begin{proposition}\label{PRws1Cat}
The first-order many-sorted rewriting systems together with the classes of first-order morphisms between first-order many-sorted rewriting systems constitute a category, denoted by $\mathsf{Rws}^{[1]}_{\mathfrak{d}}$.
\end{proposition}

\begin{proof}
That domains and codomains respect identities and compositions follows from the definitions of first-order identity morphism and first-order composition of first-order morphisms introduced in Definitions~\ref{DIdRws1} and \ref{DCompRws1}. Thus, all that remains to be proven is that the class of the first-order identity morphism at a first-order many-sorted rewriting system acts as a unit element and that the first-order composition of classes of first-order morphisms is associative.

\textsf{Unit element.}

Let $\mathbf{f}^{(1)}=(\varphi, c, (f^{(i)})_{i\in 2})$ be a first-order morphism from $\boldsymbol{\mathcal{A}}^{(1)}=(\boldsymbol{\mathcal{A}}^{(0)}, \mathcal{A}^{(1)})$ to $\boldsymbol{\mathcal{B}}^{(1)}=(\boldsymbol{\mathcal{B}}^{(0)}, \mathcal{B}^{(1)})$. We need to prove that
\begin{align*}
\left[
\mathbf{f}^{(1)}
\right] 
\circ 
\left[
\mathrm{id}^{\boldsymbol{\mathcal{A}}^{(1)}}
\right]
&=
\left[
\mathbf{f}^{(1)}
\right]
&&\mbox{ and }&
\left[
\mathrm{id}^{\boldsymbol{\mathcal{B}}^{(1)}}
\right]
\circ
\left[
\mathbf{f}^{(1)}
\right] 
&=
\left[
\mathbf{f}^{(1)}
\right].
\tag{Id}
\end{align*}

We will only prove the left hand side of (Id). The right hand side of (Id) is done similarly. According to Definition~\ref{DIdRws1}, the first-order identity morphism at $\boldsymbol{\mathcal{A}}^{(1)}$ is given by $(\mathrm{id}^{\boldsymbol{\mathcal{A}}^{(0)}}, \mathrm{ech}^{(1,\mathcal{A}^{(1)})}_{\boldsymbol{\mathcal{A}}^{(1)}})$. Moreover, according to Definition~\ref{DCompRws1} and Proposition~\ref{PCompRws1}, the first-order composition morphism is
$$
\left[
\mathbf{f}^{(1)}
\right]
\circ
\left[
\mathrm{id}^{\boldsymbol{\mathcal{A}}^{(1)}}
\right]
=
\left[
\mathbf{f}^{(1)}
\circ
\mathrm{id}^{\boldsymbol{\mathcal{A}}^{(1)}}
\right]
=
\left[
\left(
\mathbf{f}^{(0)} \circ \mathrm{id}^{\boldsymbol{\mathcal{A}}^{(0)}},
f^{(1)\flat}_{\mathrm{id}^{S}} \circ \mathrm{ech}^{(1,\mathcal{A}^{(1)})}_{\boldsymbol{\mathcal{A}}^{(1)}}
\right)
\right].
$$

According Proposition~\ref{PRws0Cat}, $\mathsf{Rws}_{\mathfrak{d}}^{(0)}$ is a category. Thus, $\mathbf{f}^{(0)} \circ \mathrm{id}^{\boldsymbol{\mathcal{A}}^{(1)}} = \mathbf{f}^{(0)}$. Moreover, the following chain of equalities holds
\allowdisplaybreaks
\begin{align*}
f^{(1)\flat}_{\mathrm{id}^{S}}
\circ
\mathrm{ech}^{(1,\mathcal{A}^{(1)})}_{\boldsymbol{\mathcal{A}}^{(1)}}
&=
f^{(1)\flat} \circ \mathrm{ech}^{(1,\mathcal{A}^{(1)})}_{\boldsymbol{\mathcal{A}}^{(1)}}
\tag{1}
\\
&=
f^{(1)}.
\tag{2}
\end{align*}

The first equality follows from the fact that, by Proposition~\ref{PMSetFunc}, $\mathrm{MSet}$ is a functor;
the second equality follows from Proposition~\ref{PPthExt}.

Thus, Equation~(Id) follows.

This proves that the class of a first-order identity morphism at a first-order many-sorted rewriting system act as the unit element for first-order composition of classes of first-order morphisms.

\textsf{Associativity.}

Let $\boldsymbol{\mathcal{A}}^{(1)}$, $\boldsymbol{\mathcal{B}}^{(1)}$, $\boldsymbol{\mathcal{C}}^{(1)}$ and $\boldsymbol{\mathcal{D}}^{(1)}$ be four first-order many-sorted rewriting systems
and let $\mathbf{f}^{(1)}$, $\mathbf{g}^{(1)}$ and $\mathbf{h}^{(1)}$ be first-order morphisms of the form
$$
\mathbf{f}^{(1)} \colon \boldsymbol{\mathcal{A}}^{(1)} \mor \boldsymbol{\mathcal{B}}^{(1)}, \,
\mathbf{g}^{(1)} \colon \boldsymbol{\mathcal{B}}^{(1)} \mor \boldsymbol{\mathcal{C}}^{(1)}
\mbox{ and }
\mathbf{h}^{(1)} \colon \boldsymbol{\mathcal{C}}^{(1)} \mor \boldsymbol{\mathcal{D}}^{(1)}
$$
where $\mathbf{f}^{(1)}$, $\mathbf{g}^{(1)}$ and $\mathbf{h}^{(1)}$ stand for
\begin{align*}
\mathbf{f}^{(1)}=(\varphi, c, (f^{(i)})_{i\in 2}), \,
\mathbf{g}^{(1)}=(\psi, d, (g^{(i)})_{i\in 2})
\\
\mbox{ and }
\mathbf{h}^{(1)}=(\theta, e, (h^{(i)})_{i\in 2}),
\mbox{ respectively.}
\end{align*}
We need to prove that 
\begin{equation}
\left[\mathbf{h}^{(1)}\right]
\circ
\left(
\left[\mathbf{g}^{(1)}\right]
\circ
\left[\mathbf{f}^{(1)}\right]
\right)
=
\left(
\left[\mathbf{h}^{(1)}\right]
\circ
\left[\mathbf{g}^{(1)}\right]
\right)
\circ
\left[\mathbf{f}^{(1)}\right].
\tag{Assoc}
\end{equation}

Following Definition~\ref{DCompRws1} and Proposition~\ref{PCompRws1}, the left hand side composition of Equation~(Assoc) is
\allowdisplaybreaks
\begin{align*}
\left[\mathbf{h}^{(1)}\right]
\circ
\left(
\left[\mathbf{g}^{(1)}\right]
\circ
\left[\mathbf{f}^{(1)}\right]
\right)
&=
\left[
\mathbf{h}^{(1)}
\circ
\left( 
\mathbf{g}^{(1)}
\circ
\mathbf{f}^{(1)}
\right)
\right]
\\
&=
\left[
\left(
\mathbf{h}^{(0)} \circ \left( \mathbf{g}^{(0)} \circ \mathbf{f}^{(0)} \right),
h^{(1)\flat}_{\psi\circ\varphi} \circ g^{(1)\flat}_{\varphi} \circ f^{(1)}
\right)
\right].
\end{align*}
Similarly, the right hand side composition of Equation~(Assoc) is
\allowdisplaybreaks
\begin{align*}
\left(
\left[\mathbf{h}^{(1)}\right]
\circ
\left[\mathbf{g}^{(1)}\right]
\right)
\circ
\left[\mathbf{f}^{(1)}\right]
&=
\left[
\left(
\mathbf{h}^{(1)}
\circ
\mathbf{g}^{(1)}
\right)
\circ
\mathbf{f}^{(1)}
\right]
\\
&=
\left[
\left(
\left(\mathbf{h}^{(0)} \circ \mathbf{g}^{(0)} \right) \circ \mathbf{f}^{(0)},
\left(h^{(1)\flat}_{\psi} \circ g^{(1)}\right)^{\flat}_{\varphi} \circ f^{(1)}
\right)
\right].
\end{align*}

According Proposition~\ref{PRws0Cat}, $\mathsf{Rws}_{\mathfrak{d}}^{(0)}$ is a category. Thus,
$$
\mathbf{h}^{(0)} \circ \left( \mathbf{g}^{(0)} \circ \mathbf{f}^{(0)} \right)
=
\left(\mathbf{h}^{(0)} \circ \mathbf{g}^{(0)} \right) \circ \mathbf{f}^{(0)}.
$$
Thus, according to Definition~\ref{DMorEqv}, all what remains to be proven is that,
$$
\mathrm{pr}^{\mathrm{Ker}(\mathrm{CH}^{(1)})}_{\boldsymbol{\mathcal{D}}^{(1)}, \theta \circ \psi \circ \varphi} 
\circ 
\left(
h^{(1)\flat}_{\psi \circ \varphi}
\circ
\left(
g^{(1)\flat}_{\varphi}
\circ
f^{(1)}
\right)
\right)^{\flat}
=
\mathrm{pr}^{\mathrm{Ker}(\mathrm{CH}^{(1)})}_{\boldsymbol{\mathcal{D}}^{(1)}, \theta \circ \psi \circ \varphi} 
\circ 
\left(
\left(
h^{(1)\flat}_{\psi}
\circ
g^{(1)}
\right)_{\varphi}^{\flat}
\circ 
f^{(1)}
\right)^{\flat}
$$

The following chain of equalities holds
\begin{flushleft}
$
\mathrm{pr}^{\mathrm{Ker}(\mathrm{CH}^{(1)})}_{\boldsymbol{\mathcal{D}}^{(1)}, \theta \circ \psi \circ \varphi} 
\circ 
\left(
h^{(1)\flat}_{\psi \circ \varphi}
\circ
\left(
g^{(1)\flat}_{\varphi}
\circ
f^{(1)}
\right)
\right)^{\flat}
$
\allowdisplaybreaks
\begin{align*}
&=
\mathrm{pr}^{\mathrm{Ker}(\mathrm{CH}^{(1)})}_{\boldsymbol{\mathcal{D}}^{(1)}, \theta \circ \psi \circ \varphi} 
\circ 
h^{(1)\flat}_{\psi \circ \varphi}
\circ
\left(
g^{(1)\flat}_{\varphi}
\circ
f^{(1)}
\right)^{\flat}
\tag{1}
\\
&=
h^{[1] @}_{\psi \circ \varphi}
\circ
\mathrm{pr}^{\mathrm{Ker}(\mathrm{CH}^{(1)})}_{\boldsymbol{\mathcal{C}}^{(1)}, \psi \circ \varphi} 
\circ
\left(
g^{(1)\flat}_{\varphi}
\circ
f^{(1)}
\right)^{\flat}
\tag{2}
\\
&=
h^{[1] @}_{\psi \circ \varphi}
\circ
\mathrm{pr}^{\mathrm{Ker}(\mathrm{CH}^{(1)})}_{\boldsymbol{\mathcal{C}}^{(1)}, \psi \circ \varphi} 
\circ
g^{(1)\flat}_{\varphi}
\circ
f^{(1)\flat}
\tag{3}
\\
&=
h^{[1] @}_{\psi \circ \varphi}
\circ
\left(
\mathrm{pr}^{\mathrm{Ker}(\mathrm{CH}^{(1)})}_{\boldsymbol{\mathcal{C}}^{(1)}, \psi} 
\circ
g^{(1)\flat}
\right)_{\varphi}
\circ
f^{(1)\flat}
\tag{4}
\\
&=
\left(
h^{[1] @}_{\psi}
\right)_{\varphi}
\circ
\left(
\mathrm{pr}^{\mathrm{Ker}(\mathrm{CH}^{(1)})}_{\boldsymbol{\mathcal{C}}^{(1)}, \psi} 
\circ
g^{(1)\flat}
\right)_{\varphi}
\circ
f^{(1)\flat}
\tag{5}
\\
&=
\left(
h^{[1] @}_{\psi}
\circ
\mathrm{pr}^{\mathrm{Ker}(\mathrm{CH}^{(1)})}_{\boldsymbol{\mathcal{C}}^{(1)}, \psi} 
\circ
g^{(1)\flat}
\right)_{\varphi}
\circ
f^{(1)\flat}
\tag{6}
\\
&=
\left(
\mathrm{pr}^{\mathrm{Ker}(\mathrm{CH}^{(1)})}_{\boldsymbol{\mathcal{D}}^{(1)}, \theta \circ \psi} 
\circ
h^{(1)\flat}_{\psi}
\circ
g^{(1)\flat}
\right)_{\varphi}
\circ
f^{(1)\flat}
\tag{7}
\\
&=
\left(
\mathrm{pr}^{\mathrm{Ker}(\mathrm{CH}^{(1)})}_{\boldsymbol{\mathcal{D}}^{(1)}, \theta \circ \psi}
\circ
\left(
h^{(1)\flat}_{\psi}
\circ
g^{(1)}
\right)^{\flat}
\right)_{\varphi}
\circ
f^{(1)\flat}
\tag{8}
\\
&=
\mathrm{pr}^{\mathrm{Ker}(\mathrm{CH}^{(1)})}_{\boldsymbol{\mathcal{D}}^{(1)}, \theta \circ \psi \circ \varphi}
\circ
\left(
h^{(1)\flat}_{\psi}
\circ
g^{(1)}
\right)^{\flat}_{\varphi}
\circ
f^{(1)\flat}
\tag{9}
\\
&=
\mathrm{pr}^{\mathrm{Ker}(\mathrm{CH}^{(1)})}_{\boldsymbol{\mathcal{D}}^{(1)}, \theta \circ \psi \circ \varphi} 
\circ 
\left(
\left(
h^{(1)\flat}_{\psi}
\circ
g^{(1)}
\right)_{\varphi}^{\flat}
\circ 
f^{(1)}
\right)^{\flat}
\tag{10}
\end{align*}
\end{flushleft}

The first equality follows from Proposition~\ref{PPthExtComp};
the second equality follows from the definition of the mapping $h^{[1]@}$, introduced in Definition~\ref{DQPthExt};
the third equality follows from Proposition~\ref{PPthExtComp};
the fourth equality follows from the fact that, according to Proposition~\ref{PDeltaPhiFunc}, $\Delta_{\varphi}$ is a convariant functor;
the fifth equality follows from the fact that, according to Proposition~\ref{PMSetFunc}, $\mathrm{MSet}$ is a contravariant functor;
the sixth equality follows from the fact that, according to Proposition~\ref{PDeltaPhiFunc}, $\Delta_{\varphi}$ is a convariant functor;
the seventh equality follows  the definition of the mapping $h^{[1]@}$, introduced in Definition~\ref{DQPthExt};
the eighth equality follows from Proposition~\ref{PPthExtComp};
the ninth equality follows from the fact that, according to Proposition~\ref{PDeltaPhiFunc}, $\Delta_{\varphi}$ is a convariant functor;
finally, the last equality follows from Proposition~\ref{PPthExtComp}.

Thus, Equation (Assoc) follows.

This proves that first-order composition of classes of first-order morphisms between is associative.

This proves that $\mathsf{Rws}^{[1]}_{\mathfrak{d}}$ is a category.
\end{proof}

\section{The category $\mathsf{Tw}_{\mathfrak{d}}^{[1]}$}

The aim of this section is to define the category of first-order towers and show that it is equivalent to the category of rewriting systems. To begin with, we show that the quotient constructions of first-order equivalent morphisms are equal.

\begin{proposition}\label{PMorEqv}
Let $\mathbf{f}^{(1)}$ and $\mathbf{f}'^{(1)}$ be first-order morphisms from $\boldsymbol{\mathcal{A}}^{(1)}$ to $\boldsymbol{\mathcal{B}}^{(1)}$. If $\mathbf{f}^{(1)}$ and $\mathbf{f}'^{(1)}$ are first-order equivalent, then 
$[\mathbf{Pth}_{\boldsymbol{\mathcal{B}}^{(1)}}^{\mathbf{f}^{(1)}}] = [\mathbf{Pth}_{\boldsymbol{\mathcal{B}}^{(1)}}^{\mathbf{f}'^{(1)}}]$ and $f^{[1] @} = f'^{[1] @}$.
\end{proposition}

\begin{proof}
Let $\mathbf{f}^{(1)} = (\mathbf{f}^{(0)}, f^{(1)})$ and $\mathbf{f}'^{(1)} = (\mathbf{f}'^{(0)}, f'^{(1)})$ be first-order morphisms from $\boldsymbol{\mathcal{A}}^{(1)}$ to $\boldsymbol{\mathcal{B}}^{(1)}$ such that $\mathbf{f}^{(1)} \cong^{(1)} \mathbf{f}'^{(1)}$. Thus, following Definition~\ref{DMorEqv}, $\mathbf{f}^{(0)} = \mathbf{f}'^{(0)}$ and $\mathrm{pr}^{\mathrm{Ker}(\mathrm{CH}^{(1)})}_{\boldsymbol{\mathcal{B}}^{(1)}, \varphi} \circ f^{(1)\flat} = \mathrm{pr}^{\mathrm{Ker}(\mathrm{CH}^{(1)})}_{\boldsymbol{\mathcal{B}}^{(1)}, \varphi} \circ f'^{(1)\flat}$. Therefore, following Definition~\ref{DQPthExt}, 
$$
f^{[1] @}
=
\left(
\mathrm{pr}^{\mathrm{Ker}(\mathrm{CH}^{(1)})}_{\boldsymbol{\mathcal{B}}^{(1)}, \varphi} \circ f^{(1)\flat}
\right)^{\natural}
=
\left(
\mathrm{pr}^{\mathrm{Ker}(\mathrm{CH}^{(1)})}_{\boldsymbol{\mathcal{B}}^{(1)}, \varphi} \circ f'^{(1)\flat}
\right)^{\natural}
=
f'^{[1] @}.
$$

Finally, we prove that the many-sorted partial $\Sigma^{\boldsymbol{\mathcal{A}}^{(1)}}$-algebras $[\mathbf{Pth}_{\boldsymbol{\mathcal{B}}^{(1)}}^{\mathbf{f}^{(1)}}]$ and $[\mathbf{Pth}_{\boldsymbol{\mathcal{B}}^{(1)}}^{\mathbf{f}'^{(1)}}]$ are equal. To this end, we consider the different items stated in Proposition~\ref{PQPthBCatAlg}.

\textsf{(1)}
The underlying $S$-sorted set of both algebras is $[\mathrm{Pth}_{\boldsymbol{\mathcal{B}}^{(1)}}]_{\varphi}$.

\textsf{(2)}
For every $(\mathbf{s}, s)\in S^{\ast}\times S$ and every operation symbol $\sigma\in\Sigma_{\mathbf{s},s}$, it follows that
$$
\sigma^{[\mathbf{Pth}_{\boldsymbol{\mathcal{B}}^{(1)}}^{\mathbf{f}^{(1)}}]}
=
\sigma^{\mathbf{c}_{\mathfrak{d}}^{\ast}([\mathbf{Pth}_{\boldsymbol{\mathcal{B}}^{(1)}}^{(0,1)}])}
=
\sigma^{[\mathbf{Pth}_{\boldsymbol{\mathcal{B}}^{(1)}}^{\mathbf{f}'^{(1)}}]}.
$$

\textsf{(3)}
For every sort $s$ in $S$ and every rewrite rule $\mathfrak{p} \in \mathcal{A}_{s}^{(1)}$, the following chain of equalities holds
\allowdisplaybreaks
\begin{align*}
\mathfrak{p}^{[\mathbf{Pth}_{\boldsymbol{\mathcal{B}}^{(1)}}^{\mathbf{f}^{(1)}}]}
&=
\mathrm{pr}^{\mathrm{Ker}(\mathrm{CH}^{(1)})}_{\boldsymbol{\mathcal{B}}^{(1)}, \varphi(s)}
\circ
f^{(1)\flat}_{(s)}
\circ
\mathrm{ech}^{(1, \mathcal{A}^{(1)})}_{\boldsymbol{\mathcal{A}}^{(1)}, s}
\left(
\mathfrak{p}
\right)
\tag{1}
\\
&=
\mathrm{pr}^{\mathrm{Ker}(\mathrm{CH}^{(1)})}_{\boldsymbol{\mathcal{B}}^{(1)}, \varphi(s)}
\circ
f'^{(1)\flat}_{(s)}
\circ
\mathrm{ech}^{(1, \mathcal{A}^{(1)})}_{\boldsymbol{\mathcal{A}}^{(1)}, s}
\left(
\mathfrak{p}
\right)
\tag{2}
\\
&=
\mathfrak{p}^{[\mathbf{Pth}_{\boldsymbol{\mathcal{B}}^{(1)}}^{\mathbf{f}'^{(1)}}]}.
\tag{3}
\end{align*}

The first equality unravels the definition of the interpretation of the constant symbol $\mathfrak{p}$ in the partial $\Sigma^{\boldsymbol{\mathcal{A}}^{(1)}}$-algebra $[\mathbf{Pth}_{\boldsymbol{\mathcal{B}}^{(1)}}^{\mathbf{f}^{(1)}}]$;
the second equality follows from the fact that $\mathbf{f}^{(1)} \cong^{(1)} \mathbf{f}'^{(1)}$;
finally, the last equality recovers the definition of the interpretation of the constant symbol $\mathfrak{p}$ in the partial $\Sigma^{\boldsymbol{\mathcal{A}}^{(1)}}$-algebra $[\mathbf{Pth}_{\boldsymbol{\mathcal{B}}^{(1)}}^{\mathbf{f}'^{(1)}}]$;

\textsf{(4)}
For every sort $s$ in $S$, it follows that
\begin{align*}
\mathrm{sc}_{s}^{0[\mathbf{Pth}_{\boldsymbol{\mathcal{B}}^{(1)}}^{\mathbf{f}^{(1)}}]}
&=
\mathrm{sc}_{\varphi(s)}^{0[\mathbf{Pth}_{\boldsymbol{\mathcal{B}}^{(1)}}]}
=
\mathrm{sc}_{s}^{0[\mathbf{Pth}_{\boldsymbol{\mathcal{B}}^{(1)}}^{\mathbf{f}'^{(1)}}]};
\\
\mathrm{tg}_{s}^{0[\mathbf{Pth}_{\boldsymbol{\mathcal{B}}^{(1)}}^{\mathbf{f}^{(1)}}]}
&=
\mathrm{tg}_{\varphi(s)}^{0[\mathbf{Pth}_{\boldsymbol{\mathcal{B}}^{(1)}}]}
=
\mathrm{tg}_{s}^{0[\mathbf{Pth}_{\boldsymbol{\mathcal{B}}^{(1)}}^{\mathbf{f}'^{(1)}}]}.
\end{align*}

\textsf{(5)}
Similarly, for every sort $s$ in $S$, it follows that
$$
\circ_{s}^{0[\mathbf{Pth}_{\boldsymbol{\mathcal{B}}^{(1)}}^{\mathbf{f}^{(1)}}]}
=
\circ_{\varphi(s)}^{0[\mathbf{Pth}_{\boldsymbol{\mathcal{B}}^{(1)}}]}
=
\circ_{s}^{0[\mathbf{Pth}_{\boldsymbol{\mathcal{B}}^{(1)}}^{\mathbf{f}'^{(1)}}]}.
$$
\end{proof}

We now define, given a first-order rewriting system, its associated first-order tower and the notion of morphism of first-order towers.

\begin{definition}\label{DTw1}
Given a first-order many-sorted rewriting system $\boldsymbol{\mathcal{A}}^{(1)} = (\boldsymbol{\mathcal{A}}^{(0)}, \mathcal{A}^{(1)})$, the \emph{first-order tower} associated to $\boldsymbol{\mathcal{A}}^{(1)}$ is 
$$
\mathbb{A}^{(1)}
=
\left(
\boldsymbol{\mathcal{A}}^{(1)}, 
\mathbf{T}_{\boldsymbol{\mathcal{E}}^{\boldsymbol{\mathcal{A}}^{(1)}}}\left(\mathbf{Pth}_{\boldsymbol{\mathcal{A}}^{(1)}}\right)
\right)
$$
where $\mathbf{T}_{\boldsymbol{\mathcal{E}}^{\boldsymbol{\mathcal{A}}^{(1)}}}(\mathbf{Pth}_{\boldsymbol{\mathcal{A}}^{(1)}})$ is the free $\Sigma^{\boldsymbol{\mathcal{A}}^{(1)}}$-algebra on the QE-variety $\mathcal{V}(\boldsymbol{\mathcal{E}}^{\boldsymbol{\mathcal{A}}^{(1)}})$. 

A \emph{morphisms of first-order towers} from $\mathbb{A}^{(1)}=(\boldsymbol{\mathcal{A}}^{(1)}, \mathbf{T}_{\boldsymbol{\mathcal{E}}^{\boldsymbol{\mathcal{A}}^{(1)}}}(\mathbf{Pth}_{\boldsymbol{\mathcal{A}}^{(1)}}))$ to $\mathbb{B}^{(1)}=(\boldsymbol{\mathcal{B}}^{(1)}, \mathbf{T}_{\boldsymbol{\mathcal{E}}^{\boldsymbol{\mathcal{B}}^{(1)}}}(\mathbf{Pth}_{\boldsymbol{\mathcal{B}}^{(1)}}))$, or simply a \emph{first-order tower morphism}, is an ordered triple $(\mathbb{A}^{(1)}, \mathbf{f}^{[1] @}, \mathbb{B}^{(1)})$, denoted by $\mathbf{f}^{[1] @} \colon \mathbb{A}^{(1)} \mor \mathbb{B}^{(1)}$ for short, in which $\mathbf{f}^{[1] @}$ is the ordered pair $([\mathbf{f}^{(1)}], f^{[1] @})$ where $[\mathbf{f}^{(1)}]$ is the first-order morphism class of a first-order morphism $\mathbf{f}^{(1)}$ from $\boldsymbol{\mathcal{A}}^{(1)}$ to $\boldsymbol{\mathcal{B}}^{(1)}$ and $f^{[1] @}$ is the $\Sigma^{\boldsymbol{\mathcal{A}}^{(1)}}$-homomorphism from $\mathbf{T}_{\boldsymbol{\mathcal{E}}^{\boldsymbol{\mathcal{A}}^{(1)}}}(\mathbf{Pth}_{\boldsymbol{\mathcal{A}}^{(1)}})$ to $\mathbf{T}_{\boldsymbol{\mathcal{E}}^{\boldsymbol{\mathcal{B}}^{(1)}}}^{\mathbf{f}^{(1)}}(\mathbf{Pth}_{\boldsymbol{\mathcal{B}}^{(1)}})$ introduced in Definition~\ref{DQPthExt}.

Note that the mapping introduced in Definition~\ref{DQPthExt}, is a $\Sigma^{\boldsymbol{\mathcal{A}}^{(1)}}$-homomorphism from $[ \mathbf{Pth}_{\boldsymbol{\mathcal{A}}^{(1)}}]$ to $[ \mathbf{Pth}_{\boldsymbol{\mathcal{B}}^{(1)}}^{\mathbf{f}^{(1)}}]$. However, taking into account Theorems~\ref{TPthFree} and \ref{TPthBFreeB}, we obtain a unique $\Sigma^{\boldsymbol{\mathcal{A}}^{(1)}}$-homomorphism, also denoted by $f^{[1] @}$ along this subsection, from $\mathbf{T}_{\boldsymbol{\mathcal{E}}^{\boldsymbol{\mathcal{A}}^{(1)}}}(\mathbf{Pth}_{\boldsymbol{\mathcal{A}}^{(1)}})$ to $\mathbf{T}_{\boldsymbol{\mathcal{E}}^{\boldsymbol{\mathcal{B}}^{(1)}}}^{\mathbf{f}^{(1)}}(\mathbf{Pth}_{\boldsymbol{\mathcal{B}}^{(1)}})$.
\end{definition}

We define below the notions of first-order identity morphism at a first-order tower and first-order composite morphism of first-order towers.

\begin{definition}
\label{DIdTw1}
Let $\mathbb{A}^{(1)}=(\boldsymbol{\mathcal{A}}^{(1)}, \mathbf{T}_{\boldsymbol{\mathcal{E}}^{\boldsymbol{\mathcal{A}}^{(1)}}}(\mathbf{Pth}_{\boldsymbol{\mathcal{A}}^{(1)}}))$ be the first-order tower associated to the first-order many-sorted rewriting system $\boldsymbol{\mathcal{A}}^{(1)}$. The ordered pair $([\mathrm{id}^{\boldsymbol{\mathcal{A}}^{(1)}}], \mathrm{id}^{\mathbf{T}_{\boldsymbol{\mathcal{E}}^{\boldsymbol{\mathcal{A}}^{(1)}}}(\mathbf{Pth}_{\boldsymbol{\mathcal{A}}^{(1)}})})$, denoted by $\mathrm{id}^{\mathbb{A}^{(1)}}$, is the first-order identity morphism at $\mathbb{A}^{(1)}$. Recall that, $[\mathrm{id}^{\boldsymbol{\mathcal{A}}^{(1)}}]$ is the first-order morphism class of the first-order identity morphism at $\boldsymbol{\mathcal{A}}^{(1)}$ defined as
$$
\mathrm{id}^{\boldsymbol{\mathcal{A}}^{(1)}}
=
\left(
\mathrm{id}^{\boldsymbol{\mathcal{A}}^{(1)}},
\mathrm{ech}^{(1,\mathcal{A}^{(1)})}_{\boldsymbol{\mathcal{A}}^{(1)}}
\right)
$$
and, according to Proposition~\ref{PQPthExtEch}, $\mathrm{id}^{\mathbf{T}_{\boldsymbol{\mathcal{E}}^{\boldsymbol{\mathcal{A}}^{(1)}}}(\mathbf{Pth}_{\boldsymbol{\mathcal{A}}^{(1)}})}$ is its quotient path extension mapping.
\end{definition}

\begin{definition}
\label{DCompTw1}
Let $\mathbf{f}^{[1] @}=([\mathbf{f}^{(1)}], f^{[1] @})$ be a first-order tower morphism from $\mathbb{A}^{(1)}$ to $\mathbb{B}^{(1)}$ and $\mathbf{g}^{[1] @}=([\mathbf{g}^{(1)}], g^{[1] @})$ a first-order tower morphism from $\mathbb{B}^{(1)}$ to $\mathbb{C}^{(1)}$ where $\mathbb{A}^{(1)}$, $\mathbb{B}^{(1)}$ and $\mathbb{C}^{(1)}$ stand for
$$
\mathbb{A}^{(1)}=(\boldsymbol{\mathcal{A}}^{(1)}, \mathbf{T}_{\boldsymbol{\mathcal{E}}^{\boldsymbol{\mathcal{A}}^{(1)}}}(\mathbf{Pth}_{\boldsymbol{\mathcal{A}}^{(1)}})), \,
\mathbb{B}^{(1)}=(\boldsymbol{\mathcal{B}}^{(1)}, \mathbf{T}_{\boldsymbol{\mathcal{E}}^{\boldsymbol{\mathcal{B}}^{(1)}}}(\mathbf{Pth}_{\boldsymbol{\mathcal{B}}^{(1)}}))
$$
$$
\mbox{ and }
\mathbb{C}^{(1)}=(\boldsymbol{\mathcal{C}}^{(1)}, \mathbf{T}_{\boldsymbol{\mathcal{E}}^{\boldsymbol{\mathcal{C}}^{(1)}}}(\mathbf{Pth}_{\boldsymbol{\mathcal{C}}^{(1)}})),
\mbox{ respectively.}
$$
The \emph{first-order composite morphism} $\mathbf{g}^{[1] @}\circ\mathbf{f}^{[1] @}$, from $\mathbb{A}^{(1)}$ to $\mathbb{C}^{(1)}$, is
$$
\mathbf{g}^{[1] @}\circ\mathbf{f}^{[1] @}=
\left(
\left[\mathbf{g}^{(1)}\right]\circ\left[\mathbf{f}^{(1)}\right], 
g^{[1] @}_{\varphi} \circ f^{[1] @}
\right).
$$
Recall that the composition $[\mathbf{g}^{(1)}] \circ [\mathbf{f}^{(1)}]$ is defined to be the first-order morphism class $[\mathbf{g}^{(1)} \circ \mathbf{f}^{(1)}]$ and, according to Corollary~\ref{CQPthExtComp}, $g^{[1] @}_{\varphi} \circ f^{[1] @}$ is its quotient path extension mapping.
\end{definition}

We show that first-order towers and first-order morphisms between them constitute a category.

\begin{proposition}\label{PTw1Cat}
The first-order towers together with the first-order morphisms between first-order towers constitute a category, denoted by $\mathsf{Tw}_{\mathfrak{d}}^{[1]}$.
\end{proposition}

\begin{proof}
That domains and codomains respect identities and compositions follows from the definitions of first-order identity morphism and first-order composition of first-order morphisms introduced in Definition~\ref{DCompTw1}. Thus, all that remains to be proven is that the first-order identity morphism at a first-order tower acts as a unit element and that the first-order composition of morphisms of first-order towers is associative.

{\sffamily Unit element.}

Let $\mathbf{f}^{[1] @}=([\mathbf{f}^{(1)}], f^{[1] @})$ be a first-order tower morphism from $\mathbb{A}^{(1)}$ to $\mathbb{B}^{(1)}$. We need to prove that 
\begin{align*}
\mathbf{f}^{[1] @} \circ \mathrm{id}^{\mathbb{A}^{(1)}} &= \mathbf{f}^{[1] @}
&&\mbox{ and }&
\mathrm{id}^{\mathbb{B}^{(1)}} \circ \mathbf{f}^{[1] @} &= \mathbf{f}^{[1] @}.
\tag{Id}
\end{align*}

We will only prove the left hand side of (Id). The right hand side of (Id) is done similarly. Let us prove that $\mathbf{f}^{[1] @} \circ \mathrm{id}^{\mathbb{A}^{(1)}} = \mathbf{f}^{[1] @}$.

The following chain of equalities holds
\allowdisplaybreaks
\begin{align*}
\mathbf{f}^{[1] @} \circ \mathrm{id}^{\mathbb{A}^{(1)}}
&=
\left( \left[\mathbf{f}^{(1)}\right], f^{[1] @} \right) \circ \left( \left[\mathrm{id}^{\boldsymbol{\mathcal{A}}^{(1)}}\right], \mathrm{id}^{\mathbf{T}_{\boldsymbol{\mathcal{E}}^{\boldsymbol{\mathcal{A}}^{(1)}}}(\mathbf{Pth}_{\boldsymbol{\mathcal{A}}^{(1)}})} \right)
\tag{1}
\\
&=
\left(
\left[\mathbf{f}^{(1)}\right] \circ \left[\mathrm{id}^{\boldsymbol{\mathcal{A}}^{(1)}}\right]
,
f^{[1] @}_{\mathrm{id}^{S}} \circ \mathrm{id}^{\mathbf{T}_{\boldsymbol{\mathcal{E}}^{\boldsymbol{\mathcal{A}}^{(1)}}}(\mathbf{Pth}_{\boldsymbol{\mathcal{A}}^{(1)}})}
\right)
\tag{2}
\\
&=
\left(
\left[\mathbf{f}^{(1)}\right]
,
f^{[1] @} \circ \mathrm{id}^{\mathbf{T}_{\boldsymbol{\mathcal{E}}^{\boldsymbol{\mathcal{A}}^{(1)}}}(\mathbf{Pth}_{\boldsymbol{\mathcal{A}}^{(1)}})}
\right)
\tag{3}
\\
&=
\left(
\left[\mathbf{f}^{(1)}\right]
,
f^{[1] @}
\right)
\tag{4}
\\
&=
\mathbf{f}^{[1]@}
\tag{5}
\end{align*}
The first equality unravels the definitions of the first-order tower morphisms $\mathbf{f}^{(1)@}$ and $\mathrm{id}^{\mathbb{A}^{(1)}}$;
the second equality unravels the definition of first-order composition of first-order tower morphisms, introduced in Definition~\ref{DCompTw1};
the third equality follows from the fact that, according to Proposition~\ref{PMSetFunc}, $\mathrm{MSet}$ is a contravariant functor;
the fourth equality follows from the fact that $f^{[1] @} \circ \mathrm{id}^{\mathbf{T}_{\boldsymbol{\mathcal{E}}^{\boldsymbol{\mathcal{A}}^{(1)}}}(\mathbf{Pth}_{\boldsymbol{\mathcal{A}}^{(1)}})}$ is equal to $f^{[1] @}$;
finally, the last equality recovers the definition of $\mathbf{f}^{[1] @}$.

Thus, Equation (Id) follows.

This proves that first-order identities at first-order towers are the unit element for first-order composition of morphisms.

{\sffamily Associativity.}

Let $\mathbb{A}^{(1)}$, $\mathbb{B}^{(1)}$, $\mathbb{C}^{(1)}$ and $\mathbb{D}^{(1)}$ be four first-order towers and let $\mathbf{f}^{[1] @}=([\mathbf{f}^{(1)}], f^{[1] @})$, $\mathbf{g}^{[1] @}=([\mathbf{g}^{(1)}], g^{[1] @})$ and $\mathbf{h}^{[1] @}=([\mathbf{h}^{(1)}], h^{[1] @})$ be first-order tower morphisms of the form
$$
\mathbf{f}^{[1] @} \colon \mathbb{A}^{(1)} \mor \mathbb{B}^{(1)}, \,
\mathbf{g}^{[1] @} \colon \mathbb{B}^{(1)} \mor \mathbb{C}^{(1)}
\mbox{ and }
\mathbf{h}^{[1] @} \colon \mathbb{C}^{(1)} \mor \mathbb{D}^{(1)}
$$
where $\mathbf{f}^{(1)}$, $\mathbf{g}^{(1)}$ and $\mathbf{h}^{(1)}$ stand for
$$
\mathbf{f}^{(1)}=(\varphi, c, (f^{(i)})_{i\in 2}), \,
\mathbf{g}^{(1)}=(\psi, d, (g^{(i)})_{i\in 2})
\mbox{ and }
\mathbf{h}^{(1)}=(\theta, e, (h^{(i)})_{i\in 2}),
\mbox{ respectively.}
$$
We need to prove that
\begin{equation}
\mathbf{h}^{[1] @} \circ \left(\mathbf{g}^{[1] @} \circ \mathbf{f}^{[1] @}\right)
=
\left(\mathbf{h}^{[1] @} \circ \mathbf{g}^{[1] @}\right) \circ \mathbf{f}^{[1] @}.
\tag{Assoc}
\end{equation}

The following chain of equalities holds
\allowdisplaybreaks
\begin{align*}
&
\mathbf{h}^{[1] @} \circ \left(\mathbf{g}^{[1] @} \circ \mathbf{f}^{[1] @}\right)
\\
&=
\left( [\mathbf{h}^{(1)}],h^{[1] @} \right) \circ \left(\left( [\mathbf{g}^{(1)}], g^{[1] @} \right) \circ \left( [\mathbf{f}^{(1)}], f^{[1] @} \right) \right)
\tag{1}
\\
&=
\left(
[\mathbf{h}^{(1)}] \circ \left([\mathbf{g}^{(1)}] \circ [\mathbf{f}^{(1)}]\right)
,
h^{[1] @}_{\psi \circ \varphi} \circ g^{[1] @}_{\varphi} \circ f^{[1] @}
\right)
\tag{2}
\\
&=
\left(
[\mathbf{h}^{(1)}] \circ \left([\mathbf{g}^{(1)}] \circ [\mathbf{f}^{(1)}]\right)
,
(h^{[1] @}_{\psi})_{\varphi} \circ g^{[1] @}_{\varphi} \circ f^{[1] @}
\right)
\tag{3}
\\
&=
\left(
\left([\mathbf{h}^{(1)}] \circ [\mathbf{g}^{(1)}]\right) \circ [\mathbf{f}^{(1)}]
,
\left(h^{[1] @}_{\psi} \circ g^{[1] @}\right)_{\varphi} \circ f^{[1] @}
\right)
\tag{4}
\\
&=
\left( \left( [\mathbf{h}^{(1)}],h^{[1] @} \right) \circ \left( [\mathbf{g}^{(1)}], g^{[1] @} \right) \right) \circ \left( [\mathbf{f}^{(1)}], f^{[1] @} \right)
\tag{5}
\\
&=
\left(\mathbf{h}^{[1] @} \circ \mathbf{g}^{[1] @}\right) \circ \mathbf{f}^{[1] @}.
\tag{6}
\end{align*}

The first equality unravels the definitions of the first-order tower morphisms $\mathbf{f}^{[1] @}$, $\mathbf{g}^{[1] @}$ and $\mathbf{h}^{[1] @}$;
the second equality unravels the definition of first-order composition of first-order tower morphisms, introduced in Definition~\ref{DCompTw1};
the third equality follows from the fact that, by Proposition~\ref{PMSetFunc}, $\mathrm{MSet}$ is a contravariant functor;
the fourth equality follows from the fact that, by Proposition~\ref{PDeltaPhiFunc}, $\Delta_{\varphi}$ is a covariant functor;
the fifth equality recovers the definition of first-order composition of first-order tower morphisms, introduced in Definition~\ref{DCompTw1};
finally, the last equality recovers the definitions of the first-order tower morphisms $\mathbf{f}^{[1] @}$, $\mathbf{g}^{[1] @}$ and $\mathbf{h}^{[1] @}$.

Thus, Equation (Assoc) follows.

This proves that first-order composition of first-order tower morphisms is associative.

This shows that $\mathsf{Tw}_{\mathfrak{d}}^{[1]}$ is a category.
\end{proof}

We will prove below that the categories of first-order many-sorted rewriting systems and first-order towers are isomorphic. In particular, we define the assignments $V^{(1)}$ and $U^{(1)}$, we prove that they are functors between the respective categories and that they are mutually inverse functors.

\begin{definition}\label{DV1}
We let $V^{(1)}$ stand for the assignment from $\mathsf{Rws}_{\mathfrak{d}}^{[1]}$ to $\mathsf{Tw}_{\mathfrak{d}}^{[1]}$ defined as follows:
\begin{enumerate}
\item
for every first-order many-sorted rewriting system $\boldsymbol{\mathcal{A}}^{(1)}$, $V^{(1)}(\boldsymbol{\mathcal{A}}^{(1)})$ is the associated first-order tower $\mathbb{A}^{(1)} = (\boldsymbol{\mathcal{A}}^{(1)}, \mathbf{T}_{\boldsymbol{\mathcal{E}}^{\boldsymbol{\mathcal{A}}^{(1)}}}(\mathbf{Pth}_{\boldsymbol{\mathcal{A}}^{(1)}}))$, and
\item
for every first-order morphism $\mathbf{f}^{(1)} \colon \boldsymbol{\mathcal{A}}^{(1)} \mor \boldsymbol{\mathcal{B}}^{(1)}$, $V^{(1)}([\mathbf{f}^{(1)}])$ is the first-order tower morphism $\mathbf{f}^{[1] @}=([\mathbf{f}^{(1)}], f^{[1] @}) \colon \mathbb{A}^{(1)} \mor \mathbb{B}^{(1)}$.
\end{enumerate}
\end{definition}

\begin{proposition}\label{PV1Fun}
The assignment $V^{(1)}$ from $\mathsf{Rws}_{\mathfrak{d}}^{[1]}$ to $\mathsf{Tw}_{\mathfrak{d}}^{[1]}$ is a covariant functor.
\end{proposition}

\begin{proof}
That $V^{(1)}$ maps objects and morphisms of $\mathsf{Rws}_{\mathfrak{d}}^{[1]}$ to objects and morphisms of $\mathsf{Tw}_{\mathfrak{d}}^{[1]}$ follows from the definition of the assignment. Therefore, all that remains to be proven is that $V^{(1)}$ preserves identities and compositions.

{\sffamily $V^{(1)}$ preserves identities.}

Let $\boldsymbol{\mathcal{A}}^{(1)}$ be an object in $\mathsf{Rws}_{\mathfrak{d}}^{[1]}$. We need to prove that
$$
V^{(1)} \left(\left[\mathrm{id}^{\boldsymbol{\mathcal{A}}^{(1)}}\right]\right)
=
\mathrm{id}^{V^{(1)}(\boldsymbol{\mathcal{A}}^{(1)})}.
$$

The following chain of equalities holds
\allowdisplaybreaks
\begin{align*}
V^{(1)} \left(\left[\mathrm{id}^{\boldsymbol{\mathcal{A}}^{(1)}}\right]\right)
&=
\left( \left[\mathrm{id}^{\boldsymbol{\mathcal{A}}^{(1)}}\right], \mathrm{ech}^{(1,\mathcal{A}^{(1)})@}_{\boldsymbol{\mathcal{A}}^{(1)}}\right)
\tag{1}
\\
&=
\left( \left[\mathrm{id}^{\boldsymbol{\mathcal{A}}^{(1)}}\right], \mathrm{id}^{\mathbf{T}_{\boldsymbol{\mathcal{E}}^{\boldsymbol{\mathcal{A}}^{(1)}}}(\mathbf{Pth}_{\boldsymbol{\mathcal{A}}^{(1)}})}\right)
\tag{2}
\\
&=
\mathrm{id}^{(\boldsymbol{\mathcal{A}}^{(1)}, \mathbf{T}_{\boldsymbol{\mathcal{E}}^{\boldsymbol{\mathcal{A}}^{(1)}}}(\mathbf{Pth}_{\boldsymbol{\mathcal{A}}^{(1)}}))}
\tag{3}
\\
&=
\mathrm{id}^{V^{(1)}(\boldsymbol{\mathcal{A}}^{2})}.
\tag{4}
\end{align*}
The first equality unravels the definition of the assignment $V^{(1)}$ on morphisms, introduced in Definition~\ref{DV1};
the second equality follows from Proposition~\ref{PQPthExtEch}; 
the third equality recovers the definition of the first-order identity tower morphism at the first-order tower $(\boldsymbol{\mathcal{A}}^{(1)}, \mathbf{T}_{\boldsymbol{\mathcal{E}}^{\boldsymbol{\mathcal{A}}^{(1)}}}(\mathbf{Pth}_{\boldsymbol{\mathcal{A}}^{(1)}}))$, introduced in Definition~\ref{DIdTw1}; 
finally, the last equality recovers the definition of the assignment $V^{(1)}$ on objects.

This proves that $V^{(1)}$ preserves identities.

{\sffamily $V^{(1)}$ Preserves compositions.}

Let $\mathbf{f}^{(1)}=(\varphi, c, (f^{(i)})_{i\in 2}) \colon \boldsymbol{\mathcal{A}}^{(1)} \mor \boldsymbol{\mathcal{B}}^{(1)}$ and $\mathbf{g}^{(1)}=(\psi, d, (g^{(i)})_{i \in 2}) \colon \boldsymbol{\mathcal{B}}^{(1)} \mor \boldsymbol{\mathcal{C}}^{(1)}$ be first-order morphisms. We need to prove that
$$
V^{(1)} \left(\left[\mathbf{g}^{(1)}\right]
\circ
\left[\mathbf{f}^{(1)}\right]\right)
=
V^{(1)} \left(\left[\mathbf{g}^{(1)}\right]\right) \circ V^{(1)} \left(\left[\mathbf{g}^{(1)}\right]\right).
$$

The following chain of equalities holds
\begin{align*}
V^{(1)} \left(\left[\mathbf{g}^{(1)}\right]
\circ
\left[\mathbf{f}^{(1)}\right]\right)
&=
V^{(1)} \left(\left[\mathbf{g}^{(1)}\circ\mathbf{f}^{(1)}\right]\right)
\tag{1}
\\
&=
\left(
\left[\mathbf{g}^{(1)}\circ\mathbf{f}^{(1)}\right],
\left(
g^{(1)\flat}_{\varphi} \circ f^{(1)}
\right)^{@}
\right)
\tag{2}
\\
&=
\left(
\left[\mathbf{g}^{(1)}\right]\circ\left[\mathbf{f}^{(1)}\right],
\left(
g^{(1)\flat}_{\varphi} \circ f^{(1)}
\right)^{@}
\right)
\tag{3}
\\
&=
\left(
\left[\mathbf{g}^{(1)}\right]\circ\left[\mathbf{f}^{(1)}\right],
g^{[1] @}_{\varphi} \circ f^{[1] @}
\right)
\tag{4}
\\
&=
\mathbf{g}^{[1] @}
\circ
\mathbf{f}^{[1] @}
\tag{5}
\\
&=
V^{(1)} \left(\left[\mathbf{g}^{(1)}\right]\right) \circ V^{(1)} \left(\left[\mathbf{f}^{(1)}\right]\right).
\tag{6}
\end{align*}

In the just stated chain of equalities, the first equality follows by Proposition~\ref{PCompRws1};
the second equality unravels the definition of the assignment $V^{(1)}$ on morphisms, introduced in Definition~\ref{DV1};
the third equality follows by Proposition~\ref{PCompRws1};
the fourth equality follows from the fact that, according to Corollary~\ref{CQPthExtComp},
$
\left(g^{(1)\flat}_{\varphi} \circ f^{(1)}\right)^{@}
=
g^{[1] @}_{\varphi} \circ f^{[1] @};
$
the fifth equality recovers the definition of the composite first-order tower morphism, introduced in Definition~\ref{DCompTw1};
finally, the last equality recovers the definition of the assignment $V^{(1)}$ on morphisms.

This proves that $V^{(1)}$ preserves compositions.

This completes the proof. 
\end{proof}

\begin{definition}\label{DU1}
We let $U^{(1)}$ stand for the assignment from $\mathsf{Tw}_{\mathfrak{d}}^{[1]}$ to $\mathsf{Rws}_{\mathfrak{d}}^{[1]}$ defined as follows:
\begin{enumerate}
\item
for every first order tower $\mathbb{A}^{(1)} = (\boldsymbol{\mathcal{A}}^{(1)}, \mathbf{T}_{\boldsymbol{\mathcal{E}}^{\boldsymbol{\mathcal{A}}^{(1)}}}(\mathbf{Pth}_{\boldsymbol{\mathcal{A}}^{(1)}}))$, $U^{(1)}(\mathbb{A}^{(1)})$ is its underlying first order many-sorted rewriting system $\boldsymbol{\mathcal{A}}^{(1)}$, and
\item
for every morphism $\mathbf{f}^{[1] @} = ([\mathbf{f}^{(1)}], f^{[1] @}) \colon \mathbb{A}^{(1)} \mor \mathbb{B}^{(1)}$, $U^{(1)}(\mathbf{f}^{[1] @})$ is the underlying first-order morphism class $[\mathbf{f}^{(1)}] \colon \boldsymbol{\mathcal{A}}^{(1)} \mor \boldsymbol{\mathcal{B}}^{(1)}$.
\end{enumerate}
\end{definition}

\begin{proposition}\label{PU1Fun}
The assignment $U^{(1)}$ from $\mathsf{Tw}_{\mathfrak{d}}^{[1]}$ to $\mathsf{Rws}_{\mathfrak{d}}^{[1]}$ is a covariant functor.
\end{proposition}

\begin{proof}
That $U^{(1)}$ maps objects and morphisms of $\mathsf{Tw}_{\mathfrak{d}}^{[1]}$ to objects and morphisms of $\mathsf{Rws}_{\mathfrak{d}}^{[1]}$ follows from the definition of the assignment. Therefore, all that remains to be proven is that $U^{(1)}$ preserves identities and compositions.

{\sffamily $U^{(1)}$ preserves identities.}

Let $\mathbb{A}^{(1)} = (\boldsymbol{\mathcal{A}}^{(1)}, \mathbf{T}_{\boldsymbol{\mathcal{E}}^{\boldsymbol{\mathcal{A}}^{(1)}}}(\mathbf{Pth}_{\boldsymbol{\mathcal{A}}^{(1)}}))$ be an object in $\mathsf{Tw}_{\mathfrak{d}}^{[1]}$. We need to prove that
$$
U^{(1)} \left(\mathrm{id}^{\mathbb{A}^{(1)}}\right)
=
\left[\mathrm{id}^{U^{(1)}(\mathbb{A}^{(1)})}\right].
$$

The following chain of equalities holds
\begin{align*}
U^{(1)} \left(\mathrm{id}^{\mathbb{A}^{(1)}}\right)
&=
U^{(1)} \left(
\left[\mathrm{id}^{\boldsymbol{\mathcal{A}}^{(1)}}\right], 
\mathrm{id}^{\mathbf{T}_{\boldsymbol{\mathcal{E}}^{\boldsymbol{\mathcal{A}}^{(1)}}}(\mathbf{Pth}_{\boldsymbol{\mathcal{A}}^{(1)}})}
\right)
\tag{1}
\\
&=
\left[\mathrm{id}^{\boldsymbol{\mathcal{A}}^{(1)}}\right]
\tag{2}
\\
&=
\left[\mathrm{id}^{U^{(1)}(\mathbb{A}^{(1)})}\right].
\tag{3}
\end{align*}
The first equality unravels the definition of the identity first-order morphism at $\mathbb{A}^{(1)}$;
the second equality unravels the definition of the assignment $U^{(1)}$ on morphisms, introduced in Definition~\ref{DU1};
finally, the last equality recovers the definition of the assignment $U^{(1)}$ on objects.

This proves that $U^{(1)}$ preserves identities.

{\sffamily $U^{(1)}$ preserves compositions.}

Let $\mathbf{f}^{[1] @} = ([\mathbf{f}^{(1)}], f^{[1] @}) \colon \mathbb{A}^{(1)} \mor \mathbb{B}^{(1)}$ and $\mathbf{g}^{[1] @} = ([\mathbf{g}^{(1)}], g^{[1] @}) \colon \mathbb{B}^{(1)} \mor \mathbb{C}^{(1)}$ be morphisms in $\mathsf{Tw}_{\mathfrak{d}}^{[1]}$. We need to prove that
$$
U^{(1)} \left(\mathbf{g}^{[1] @} \circ \mathbf{f}^{[1] @}\right)
=
U^{(1)} \left(\mathbf{g}^{[1] @}\right) \circ U^{(1)} \left(\mathbf{f}^{[1] @}\right).
$$

The following chain of equalities holds
\begin{align*}
U^{(1)} \left(\mathbf{g}^{[1] @} \circ \mathbf{f}^{[1] @}\right)
&=
U^{(1)} \left(
\left[\mathbf{g}^{(1)}\right]\circ\left[\mathbf{f}^{(1)}\right],
g^{[1] @}_{\varphi} \circ f^{[1] @}
\right)
\tag{1}
\\
&=
\left[\mathbf{g}^{(1)}\right]\circ\left[\mathbf{f}^{(1)}\right]
\tag{2}
\\
&=
U^{(1)} \left(\mathbf{g}^{[1] @}\right) \circ U^{(1)} \left(\mathbf{f}^{[1] @}\right).
\tag{3}
\end{align*}
The first equality unravels the definition of the composition of second-order tower morphisms, introduced in Definition~\ref{DCompTw1};
the second equality unravels the definition of the assignment $U^{(1)}$ on morphisms, introduced in Definition~\ref{DU1};
finally, the last equality recovers the definition of the assignment $U^{(1)}$ on morphisms.

This proves that $U^{(1)}$ preserves compositions.

This completes the proof.
\end{proof}

\begin{proposition}\label{PV1U1Comp}
$U^{(1)} \circ V^{(1)} = \mathrm{Id}^{\mathsf{Rws}_{\mathfrak{d}}^{[1]}}$ and $V^{(1)} \circ U^{(1)} = \mathrm{Id}^{\mathsf{Tw}_{\mathfrak{d}}^{[1]}}$
\end{proposition}

\begin{proof}
For every object $\boldsymbol{\mathcal{A}}^{(1)}$ in $\mathsf{Rws}_{\mathfrak{d}}^{[1]}$, the following chain of equalities holds
\begin{align*}
U^{(1)} \left(V^{(1)} \left(\boldsymbol{\mathcal{A}}^{(1)}\right)\right)
&=
U^{(1)} \left(\boldsymbol{\mathcal{A}}^{(1)}, \mathbf{T}_{\boldsymbol{\mathcal{E}}^{\boldsymbol{\mathcal{A}}^{(1)}}}(\mathbf{Pth}_{\boldsymbol{\mathcal{A}}^{(1)}})\right)
\tag{1}
\\
&=
\boldsymbol{\mathcal{A}}^{(1)}.
\tag{2}
\end{align*}
The first equality unravels the definition of the functor $V^{(1)}$ on objects, introduced in Definition~\ref{DV1}; 
the second equality unravels the definition of the functor $U^{(1)}$ on objects, introduced in Definition~\ref{DU1}.

For every first-order morphism $\mathbf{f}^{(1)} \colon \boldsymbol{\mathcal{A}}^{(1)} \mor \boldsymbol{\mathcal{B}}^{(1)}$, the following chain of equalities holds
\begin{align*}
U^{(1)} \left(V^{(1)} \left(\left[\mathbf{f}^{(1)}\right]\right)\right)
&=
U^{(1)} \left(\left[\mathbf{f}^{(1)}\right], f^{[1] @}\right)
\tag{1}
\\
&=
\left[\mathbf{f}^{(1)}\right].
\tag{2}
\end{align*}
The first equality unravels the definition of the functor $V^{(1)}$ on morphisms, introduced in Definition~\ref{DV1};
the second equality unravels the definition of the functor $U^{(1)}$ on morphisms, introduced in Definition~\ref{DU1}.

For every object $\mathbb{A}^{(1)} = (\boldsymbol{\mathcal{A}}^{(1)}, \mathbf{T}_{\boldsymbol{\mathcal{E}}^{\boldsymbol{\mathcal{A}}^{(1)}}}(\mathbf{Pth}_{\boldsymbol{\mathcal{A}}^{(1)}}))$ in $\mathsf{Tw}_{\mathfrak{d}}^{[1]}$, the following chain of equalities holds
\begin{align*}
V^{(1)} \left(U^{(1)} \left(\mathbb{A}^{(1)}\right)\right)
&=
V^{(1)} \left(\boldsymbol{\mathcal{A}}^{(1)}\right)
\tag{1}
\\
&=
(\boldsymbol{\mathcal{A}}^{(1)}, \mathbf{T}_{\boldsymbol{\mathcal{E}}^{\boldsymbol{\mathcal{A}}^{(1)}}}(\mathbf{Pth}_{\boldsymbol{\mathcal{A}}^{(1)}}))
\tag{2}
\\
&=
\mathbb{A}^{(1)}.
\tag{3}
\end{align*}
The first equality unravels the definition of the functor $U^{(1)}$ on objects, introduced in Definition~\ref{DU1};
the second equality unravels the definition of the functor $V^{(1)}$ on objects, introduced in Definition~\ref{DV1};
finally, the last equality recovers the definition $\mathbb{A}^{(1)}$.

For every morphism $\mathbf{f}^{[1] @} = \left(\left[\mathbf{f}^{(1)}\right], f^{[1] @}\right) \colon \mathbb{A}^{(1)} \mor \mathbb{B}^{(1)}$ in $\mathsf{Tw}_{\mathfrak{d}}^{[1]}$, the following chain of equalities holds
\begin{align*}
V^{(1)} \left(U^{(1)} \left(\mathbf{f}^{[1] @}\right)\right)
&=
V^{(1)} \left(\left[\mathbf{f}^{(1)}\right]\right)
\tag{1}
\\
&=
\left(\left[\mathbf{f}^{(1)}\right], f^{[1] @}\right)
\tag{2}
\\
&=
\mathbf{f}^{[1] @}.
\tag{3}
\end{align*}
The first equality unravels the definition of the functor $U^{(1)}$ on morphisms, introduced in Definition~\ref{DU1};
the second equality unravels the definition of the functor $V^{(1)}$ on morphisms, introduced in Definition~\ref{DV1};
finally, the last equality recovers the definition $\mathbf{f}^{[1] @}$.

This completes the proof.
\end{proof}

\chapter{Morphisms of second-order many-sorted rewriting systems}\label{S3E}

Let us recall that a second-order many-sorted rewriting system is an ordered pair $\boldsymbol{\mathcal{A}}^{(2)}=(\boldsymbol{\mathcal{A}}^{(1)}, \mathcal{A}^{(2)})$, where $\boldsymbol{\mathcal{A}}^{(1)} = (\boldsymbol{\mathcal{A}}^{(0)}, \mathcal{A}^{(1)})$ is a first-order $S$-sorted rewriting system, see Definition~\ref{DRewSys}, and, for the many sorted signature $\boldsymbol{\Sigma}^{\boldsymbol{\mathcal{A}}^{(1)}} = (S, \Sigma^{\boldsymbol{\mathcal{A}}^{(1)}})$, $\mathcal{A}^{(2)}$ is a subset of
\begin{multline*}
\mathrm{Rwr}(\mathbf{\Sigma}^{\boldsymbol{\mathcal{A}}^{(1)}}, X)
=
\left(
\left\lbrace
([M]_{s}, [N]_{s}) \in [\mathrm{PT}_{\boldsymbol{\mathcal{A}}^{(1)}}]^{2}_{s} 
\middle|
\right.\right.
\\
\left.\left.
\left(
\mathrm{ip}^{(1,X)@}_{s}(M),
\mathrm{ip}^{(1,X)@}_{s}(N)
\right)
\in \mathrm{Ker}\left(\mathrm{sc}^{(0,1)}\right)_{s}
\cap\mathrm{Ker}\left(\mathrm{tg}^{(0,1)}\right)_{s}
\right\rbrace
\right)_{s\in S},
\end{multline*}
the $S$-sorted set of the \emph{second-order rewrite rules with variables in $X$}. We will also call $\boldsymbol{\mathcal{A}}^{(2)}$ a second-order $S$-sorted rewriting system. For further details see Definition~\ref{DDRewSys}.

In this chapter, given two second-order many-sorted rewriting systems $\boldsymbol{\mathcal{A}}^{(2)}=(\boldsymbol{\mathcal{A}}^{(1)}, \mathcal{A}^{(2)})$ and  $\boldsymbol{\mathcal{B}}^{(2)}=(\boldsymbol{\mathcal{B}}^{(1)}, \mathcal{B}^{(2)})$, we define the notion of second-order morphism $\mathbf{f}^{(2)}$ from  $\boldsymbol{\mathcal{A}}^{(2)}$ to  $\boldsymbol{\mathcal{B}}^{(2)}$, where $\mathbf{f}^{(2)}=(\mathbf{f}^{(1)},f^{(2)})$, $\mathbf{f}^{(1)}=(\varphi,c,(f^{(i)})_{i \in 2})$ is a first-order morphism from $\boldsymbol{\mathcal{A}}^{(1)}$ to $\boldsymbol{\mathcal{B}}^{(1)}$, see Definition~\ref{DRws1Mor}, and $f^{(2)}$ is an $S$-sorted mapping from $\mathcal{A}^{(2)}$ to $\mathrm{Pth}_{\boldsymbol{\mathcal{B}}^{(2)},\varphi}$ with some compatibility properties with sources and targets.
We next define a structure of $\Sigma$-algebra on $\mathrm{Pth}_{\boldsymbol{\mathcal{B}}^{(2)}, \varphi}$ and $\llbracket \mathrm{Pth}_{\boldsymbol{\mathcal{B}}^{(2)}, \varphi} \rrbracket$, which we will denote by $\mathbf{Pth}_{\boldsymbol{\mathcal{B}}^{(2)}}^{\mathbf{f}^{(2)}(0,2)}$ and $\llbracket \mathbf{Pth}_{\boldsymbol{\mathcal{B}}^{(2)}}^{\mathbf{f}^{(2)}(0,2)} \rrbracket$, respectively.
Then, given a second-order morphism $\mathbf{f}^{(2)}$ from $\boldsymbol{\mathcal{A}}^{(2)}$ to $\boldsymbol{\mathcal{B}}^{(2)}$, we consider an extension of the mapping $f^{(2)}$ to the $S$-sorted set $\mathrm{Pth}_{\boldsymbol{\mathcal{A}}^{(2)}}$, which we have called the second-order path extension mapping, and denote it by $f^{(2)\flat}$.
After that, we prove that the path extension mapping is a $\Sigma$-homomorphism.
Moreover, we show its relation with the $S$-sorted mappings $\mathrm{sc}^{(0,2)}$, $\mathrm{tg}^{(0,2)}$ and $\mathrm{ip}^{(2,0)\sharp}$.
Finally, we define a structure of partial $\Sigma^{\boldsymbol{\mathcal{A}}^{(1)}}$-algebra on $\mathrm{Pth}_{\boldsymbol{\mathcal{B}}^{(2)}, \varphi}$ and $\llbracket \mathrm{Pth}_{\boldsymbol{\mathcal{B}}^{(2)}, \varphi} \rrbracket$, which we will denote by $\mathbf{Pth}_{\boldsymbol{\mathcal{B}}^{(2)}}^{\mathbf{f}^{(2)}(0,2)}$ and $\llbracket \mathbf{Pth}_{\boldsymbol{\mathcal{B}}^{(2)}}^{\mathbf{f}^{(2)}(0,2)} \rrbracket$, respectively and show that the second-order path extension mapping is a $\Sigma^{\boldsymbol{\mathcal{A}}^{(1)}}$-homomorphism.

\begin{definition}\label{DRws2Mor}
A \emph{morphism} of second-order many-sorted rewriting systems, or simply a \emph{second-order morphism}, from $\boldsymbol{\mathcal{A}}^{(2)}=(\boldsymbol{\mathcal{A}}^{(1)}, \mathcal{A}^{(2)})$ to $\boldsymbol{\mathcal{B}}^{(2)}=(\boldsymbol{\mathcal{B}}^{(1)}, \mathcal{B}^{(2)})$ is an ordered triple $(\boldsymbol{\mathcal{A}}^{(2)}, \mathbf{f}^{(2)}, \boldsymbol{\mathcal{B}}^{(2)})$, denoted by $\mathbf{f}^{(2)} \colon \boldsymbol{\mathcal{A}}^{(2)} \mor \boldsymbol{\mathcal{B}}^{(2)}$ for short, in which $\mathbf{f}^{(2)}=(\varphi, c, (f^{(i)})_{i\in 3})$ is the ordered pair where 
\begin{enumerate}
\item
$\mathbf{f}^{(1)} = (\varphi, c, (f^{(i)})_{i\in 2})$, the \emph{underlying first-order morphism} of $\mathbf{f}^{(2)}$, is a first-order morphism from $\boldsymbol{\mathcal{A}}^{(1)}$ to $\boldsymbol{\mathcal{B}}^{(1)}$, introduced in Definition~\ref{DRws1Mor}.
$$
\mathbf{f}^{(1)} \colon
\boldsymbol{\mathcal{A}}^{(1)}
\mor
\boldsymbol{\mathcal{B}}^{(1)}
$$
\item 
$f^{(2)} \colon \mathcal{A}^{(2)} \mor \mathrm{Pth}_{\boldsymbol{\mathcal{B}}^{(2)},\varphi}$ is an $S$-sorted mapping, where $\mathrm{Pth}_{\boldsymbol{\mathcal{B}}^{(2)},\varphi}$ is the $S$-sorted set $(\mathrm{Pth}_{\boldsymbol{\mathcal{B}}^{(2)},\varphi(s)})_{s\in S}$, satisfying that, for every $s\in S$ and every rewrite rule $\mathfrak{p}^{(2)}=([M]_{s},[N]_{s})\in\mathcal{A}^{(2)}_{s}$, we have that
$$
f_{s}^{(2)}\left(\mathfrak{p}^{(2)}\right)
\in
\mathrm{Pth}_{\boldsymbol{\mathcal{B}}^{(2)}, \varphi(s)}\left(
f^{[1]@}_{s}\left(
\left[M\right]_{s}
\right), 
f^{[1]@}_{s}\left(
\left[N\right]_{s}
\right)
\right),
$$
or, what is equivalent,
\begin{align*}
\mathrm{sc}^{([1],2)}_{\boldsymbol{\mathcal{B}}^{(2)},\varphi(s)}\left(
f_{s}^{(2)}\left(
\mathfrak{p}^{(2)}
\right)
\right)
&=
f^{[1]@}_{s}\left(
\left[M\right]_{s}
\right)
\mbox{ and}
\\
\mathrm{tg}^{([1],2)}_{\boldsymbol{\mathcal{B}}^{(2)},\varphi(s)}\left(
f_{s}^{(2)}\left(
\mathfrak{p}^{(2)}
\right)
\right)
&=
f^{[1]@}_{s}\left(
\left[N\right]_{s}
\right)
\end{align*}
\end{enumerate}

The alternative notation $\mathbf{f}^{(2)} = (\mathbf{f}^{(1)}, f^{(2)})$ will also be used.

\end{definition}

\begin{remark}
\label{RDDSigmaAlg}
Let $\mathbf{f}^{(2)}=(\varphi, c, (f^{(i)})_{i\in 3})$ be a second-order morphism from $\boldsymbol{\mathcal{A}}^{(2)}$ to $\boldsymbol{\mathcal{B}}^{(2)}$.
Then the $S$-sorted sets $\mathrm{Pth}_{\boldsymbol{\mathcal{B}}^{(2)}, \varphi}$ and $\llbracket\mathrm{Pth}_{\boldsymbol{\mathcal{B}}^{(2)}}\rrbracket_{\varphi}$ are equipped, in a natural way, with a structure of $\Sigma$-algebra, namely, $\mathbf{c}_{\mathfrak{d}}^{\ast}(\mathbf{Pth}_{\boldsymbol{\mathcal{B}}^{(2)}}^{(0,2)})$ and $\mathbf{c}_{\mathfrak{d}}^{\ast}([\mathbf{Pth}_{\boldsymbol{\mathcal{B}}^{(2)}}^{(0,2)}])$.
In order to unify the notation, we will denote these $\Sigma$-algebras by $\mathbf{Pth}_{\boldsymbol{\mathcal{B}}^{(2)}}^{\mathbf{f}^{(2)}(0,2)}$ and $\llbracket\mathbf{Pth}_{\boldsymbol{\mathcal{B}}^{(2)}}^{\mathbf{f}^{(2)}(0,2)}\rrbracket$, respectively.
Moreover, the $S$-sorted mapping $\mathrm{pr}_{\boldsymbol{\mathcal{B}}^{(2)}, \varphi}^{\llbracket\cdot\rrbracket} = \mathbf{c}_{\mathfrak{d}}^{\ast}(\mathrm{pr}_{\boldsymbol{\mathcal{B}}^{(2)}}^{\llbracket\cdot\rrbracket})$ is a $\Sigma$-homomorphism from  $\mathbf{Pth}_{\boldsymbol{\mathcal{B}}^{(2)}}^{\mathbf{f}^{(2)}(0,2)}$ to $\llbracket\mathbf{Pth}_{\boldsymbol{\mathcal{B}}^{(2)}}^{\mathbf{f}^{(2)}(0,2)}\rrbracket$ by Propositions~\ref{PDVDZPr} and \ref{PFunSig}.
Similarly, $\mathrm{sc}^{(0,2)}_{\boldsymbol{\mathcal{B}}^{(2)}, \varphi} = \mathbf{c}_{\mathfrak{d}}^{\ast}(\mathrm{sc}^{(0,2)}_{\boldsymbol{\mathcal{B}}^{(2)}})$ and $\mathrm{tg}^{(0,2)}_{\boldsymbol{\mathcal{B}}^{(2)}, \varphi} = \mathbf{c}_{\mathfrak{d}}^{\ast}(\mathrm{tg}^{(0,2)}_{\boldsymbol{\mathcal{B}}^{(2)}})$ are $\Sigma$-homomorphisms from $\mathbf{Pth}_{\boldsymbol{\mathcal{B}}^{(2)}}^{\mathbf{f}^{(2)}(0,2)}$ to $\mathbf{c}_{\mathfrak{d}}^{\ast}(\mathbf{T}_{\Lambda}(Y))$ by Propositions~\ref{PDZHom} and \ref{PFunSig}.
Finally, $\mathrm{ip}^{(2,0)\sharp}_{\boldsymbol{\mathcal{B}}^{(2)}, \varphi} = \mathbf{c}_{\mathfrak{d}}^{\ast} (\mathrm{ip}^{(2,0)\sharp}_{\boldsymbol{\mathcal{B}}^{(2)}})$ is a  $\Sigma$-homomorphism from $\mathbf{c}_{\mathfrak{d}}^{\ast}(\mathbf{T}_{\Lambda}(Y))$ to $\mathbf{Pth}_{\boldsymbol{\mathcal{B}}^{(2)}}^{\mathbf{f}^{(2)}(0,2)}$ by Propositions~\ref{PDZHomIp} and \ref{PFunSig}.
\end{remark}

\section{Second-order path extension mapping of a morphism}

In this section we prove that, for a second-order morphism $\mathbf{f}^{(2)}$, there exists an extension of $f^{(2)}$, the \emph{second-order path extension mapping}, to the set of second-order paths $\mathrm{Pth}_{\boldsymbol{\mathcal{A}}^{(2)}}$.

\begin{proposition}
\label{PDPthExt}
Let $\mathbf{f}^{(2)}=(\varphi, c, (f^{(i)})_{i\in 3})$ be a second-order morphism from $\boldsymbol{\mathcal{A}}^{(2)}$ to $\boldsymbol{\mathcal{B}}^{(2)}$. Then there exists an $S$-sorted mapping $f^{(2)\flat}$ from $\mathrm{Pth}_{\boldsymbol{\mathcal{A}}^{(2)}}$ to $\mathrm{Pth}_{\boldsymbol{\mathcal{B}}^{(2)}, \varphi}$, which we call the \emph{second-order path extension mapping of $f^{(2)}$}, satisfying that
\begin{multicols}{2}
\begin{enumerate}
\item $\mathrm{sc}^{([1],2)}_{\boldsymbol{\mathcal{B}}^{(2)},\varphi}\circ f^{(2)\flat}=f^{[1]@}\circ \mathrm{sc}^{([1],2)}_{\boldsymbol{\mathcal{A}}^{(2)}}$;
\item $\mathrm{tg}^{([1],2)}_{\boldsymbol{\mathcal{B}}^{(2)},\varphi}\circ f^{(2)\flat}=f^{[1]@}\circ \mathrm{tg}^{([1],2)}_{\boldsymbol{\mathcal{A}}^{(1)}}$;
\item
$f^{(2)\flat} \circ \mathrm{ip}_{\boldsymbol{\mathcal{A}}^{(2)}}^{(2, [1])\sharp} = \mathrm{ip}_{\boldsymbol{\mathcal{B}}^{(2)}, \varphi}^{(2, [1])\sharp} \circ f^{[1]@}$;
\item
$f^{(2)\flat} \circ \mathrm{ech}^{(2,\mathcal{A}^{(2)})}_{\boldsymbol{\mathcal{A}}^{(2)}} = f^{(2)}$.
\end{enumerate}
\end{multicols}
\end{proposition}

\begin{proof}
Let us define the $S$-sorted mapping $f^{(2)\flat}$ by Artinian recursion on $(\coprod \mathrm{Pth}_{\boldsymbol{\mathcal{A}}^{(2)}}, \leq_{\mathbf{Pth}_{\boldsymbol{\mathcal{A}}^{(2)}}})$ as follows.

{\sffamily Base step of the Artinian recursion.}

Let $(\mathfrak{P}^{(2)}, s)$ be a minimal element of $(\coprod \mathrm{Pth}_{\boldsymbol{\mathcal{A}}^{(2)}}, \leq_{\mathbf{Pth}_{\boldsymbol{\mathcal{A}}^{(2)}}})$. Then, by Proposition~\ref{PDMinimal}, the path $\mathfrak{P}^{(2)}$ is either~(1) an $(2,[1])$-identity second-order path or~(2) a second-order echelon.

If~(1), i.e., if $\mathfrak{P}^{(2)}$ is an $(2,[1])$-identity second-order path, then $\mathfrak{P}^{(2)}=\mathrm{ip}^{(2,[1])\sharp}_{\boldsymbol{\mathcal{A}}^{(2)}, s}([P]_{s})$ for some term $[P]_{s} \in [\mathrm{PT}_{\boldsymbol{\mathcal{A}}^{(2)}}]_{s}$. We define $f^{(2)\flat}_{s}(\mathfrak{P}^{(2)})$ to be the $(2,[1])$-identity second-order path at $f_{s}^{[1]@}([P]_{s})$ which is a path term class in $[\mathrm{PT}_{\boldsymbol{\mathcal{B}}^{(1)}}(Y)]_{\varphi(s)}$, i.e.,
$$
f^{(2)\flat}_{s}(\mathfrak{P}^{(2)}) = \mathrm{ip}^{(2,[1])\sharp}_{\boldsymbol{\mathcal{B}}^{(2)}, \varphi(s)}\left(
f_{s}^{[1]@}\left(
[P]_{s}
\right)
\right).
$$

If~(2), i.e., if $\mathfrak{P}^{(2)}$ is an echelon associated to a rewrite rule $\mathfrak{p}^{(2)}=([M]_{s}, [N]_{s})$, that is, if $\mathfrak{P}^{(2)}$ has the form
$$
\xymatrix@C=105pt{
\mathfrak{P}^{(2)}: [M]_{s}
\ar@{=>}[r]^-{\text{\Small{$\left(\mathfrak{p}^{(2)},\mathrm{id}^{\mathrm{T}_{\Sigma^{\boldsymbol{\mathcal{A}}^{(1)}}}(X)_{s}}\right)$}}}
&
[N]_{s}
},
$$
then we define $f_{s}^{(2)\flat}\left(\mathfrak{P}^{(2)}\right)$ to be the image of $\mathfrak{p}^{(2)}$, the unique rewrite rule appearing in $\mathfrak{P}^{(2)}$, under the mapping $f_{s}^{(2)}$, i.e.,
$$
f^{(2)\flat}_{s}\left(\mathfrak{P}^{(2)}\right)
=
f_{s}^{(2)}\left(\mathfrak{p}^{(2)}\right).
$$

This completes the base step of the Artinian recursion.

{\sffamily Inductive step of the Artinian recursion.}

Let $(\mathfrak{P}^{(2)},s)$ be a non-minimal element of $(\coprod\mathrm{Pth}_{\boldsymbol{\mathcal{A}}^{(2)}}, \leq_{\mathbf{Pth}_{\boldsymbol{\mathcal{A}}^{(2)}}})$. We can assume that $\mathfrak{P}^{(2)}$ is a not a $(2,[1])$-identity second-order path, since those paths already have an image for the mapping $f^{(2)\flat}$. Let us suppose that, for every sort $t\in S$ and every second-order path $\mathfrak{Q}^{(2)}\in\mathrm{Pth}_{\boldsymbol{\mathcal{A}}^{(1)},t}$, if $(\mathfrak{Q}^{(2)},t)<_{\mathbf{Pth}_{\boldsymbol{\mathcal{A}}^{(1)}}}(\mathfrak{P}^{(2)},s)$, then the value of the mapping $f^{(2)\flat}$ at $\mathfrak{Q}^{(2)}$, i.e., $f^{(2)\flat}_{t}(\mathfrak{Q}^{(2)})$, has already been defined. Moreover, assume that, for every $t \in S$ and every second-order path $\mathfrak{Q}^{(2)}\in\mathrm{Pth}_{\boldsymbol{\mathcal{A}}^{(2)}, t}$, the definition of $f^{(2)\flat}_{t}(\mathfrak{Q}^{(2)})$ satisfies the following equality
\allowdisplaybreaks
\begin{align*}
\mathrm{sc}_{\boldsymbol{\mathcal{B}}^{(2)},\varphi(t)}^{([1],2)}\left(f^{(2)\flat}_{t}\left(\mathfrak{Q}^{(2)}\right)\right) 
&=
f^{[1]@}_{t}\left(\mathrm{sc}_{\boldsymbol{\mathcal{A}}^{(2)},t}^{([1],2)}\left(\mathfrak{Q}^{(2)}\right)\right);
\\
\mathrm{tg}_{\boldsymbol{\mathcal{B}}^{(2)},\varphi(t)}^{([1],2)}\left(f^{(2)\flat}_{t}\left(\mathfrak{Q}^{(2)}\right)\right)
&=
f^{[1]@}_{t}\left(\mathrm{tg}_{\boldsymbol{\mathcal{A}}^{(1)},t}^{([1],2)}\left(\mathfrak{Q}^{(2)}\right)\right).
\end{align*}

By Lemma~\ref{LDOrdI}, we have that $\mathfrak{P}^{(2)}$ is either~(1) a path of length strictly greater than one containing at least one second-order echelon or~(2) an echelonless second-order path.

If~(1), i.e., if $\mathfrak{P}^{(2)}$ is a path of length strictly greater than one containing at least one second-order echelon, then let $i\in \bb{\mathfrak{P}^{(2)}}$ be the first index for which the one-step subpath $\mathfrak{P}^{(2),i,i}$ of $\mathfrak{P}^{(2)}$ is a second-order echelon. We consider different cases for $i$ according to the cases presented in Definition~\ref{DDOrd}.

If $i=0$, we have that the pairs $(\mathfrak{P}^{(2),0,0},s)$ and $(\mathfrak{P}^{(2),1,\bb{\mathfrak{P}^{(2)}}-1},s)$ $\prec_{\mathbf{Pth}_{\boldsymbol{\mathcal{A}}^{(2)}}}$-precede the pair $(\mathfrak{P}^{(2)},s)$. Therefore, the values of the mapping $f^{(2)\flat}$ at $\mathfrak{P}^{(2),0,0}$ and $\mathfrak{P}^{(2),1,\bb{\mathfrak{P}^{(2)}}-1}$, respectively, have already been defined. In particular, the following chain of equalities holds
\allowdisplaybreaks
\begin{align*}
\mathrm{sc}^{([1],2)}_{\boldsymbol{\mathcal{B}}^{(2)},\varphi(s)}\left(
f_{s}^{(2)\flat}\left(
\mathfrak{P}^{(2),1,\bb{\mathfrak{P}}-1}
\right)
\right)
&=
f^{[1]@}_{s}\left(
\mathrm{sc}^{([1],2)}_{\boldsymbol{\mathcal{A}}^{(2)},s}\left(
\mathfrak{P}^{(2),1,\bb{\mathfrak{P}}-1}
\right)
\right)
\tag{1}
\\
&=
f^{[1]@}_{s}\left(
\mathrm{tg}^{([1],2)}_{\boldsymbol{\mathcal{A}}^{(2)},s}\left(
\mathfrak{P}^{(2),0,0}
\right)
\right)
\tag{2}
\\
&=
\mathrm{tg}^{([1],2)}_{\boldsymbol{\mathcal{B}}^{(2)},\varphi(s)}\left(
f_{s}^{(2)\flat}\left(
\mathfrak{P}^{(2),0,0}
\right)
\right).
\tag{3}
\end{align*}

The first equality follows from the assumption of the Artinian recursion;
the second equality follows from Proposition~\ref{PDPthRecons};
finally, the last equality follows from the assumption of the Artinian recursion.

In this case, we set
\allowdisplaybreaks
\begin{align*}
f^{(2)\flat}_{s}\left(
\mathfrak{P}^{(2)}
\right)
&=
f^{(2)\flat}_{s}\left(
\mathfrak{P}^{(2),1,\bb{\mathfrak{P}^{(2)}}-1}
\right)
\circ_{\varphi(s)}^{1\mathbf{Pth}_{\boldsymbol{\mathcal{B}}^{(2)}}}
f^{(2)\flat}_{s}\left(
\mathfrak{P}^{(2),0,0}
\right).
\end{align*}
The above chain of equalities justifies that such composite exists.

If $i\neq 0$, we have that the pairs $(\mathfrak{P}^{(2),0,i-1},s)$ and $(\mathfrak{P}^{(2),i,\bb{\mathfrak{P}^{(2)}}-1},s)$ $\prec_{\mathbf{Pth}_{\boldsymbol{\mathcal{A}}^{(2)}}}$-precede the pair $(\mathfrak{P}^{(2)},s)$. Therefore, the values of the mapping $f^{(2)\flat}$ at $\mathfrak{P}^{(2),0,i-1}$ and $\mathfrak{P}^{(2),i,\bb{\mathfrak{P}^{(2)}}-1}$, respectively, have already been defined. In particular, the following chain of equalities holds
\allowdisplaybreaks
\begin{align*}
\mathrm{sc}^{([1],2)}_{\boldsymbol{\mathcal{B}}^{(2)},\varphi(s)}\left(
f_{s}^{(2)\flat}\left(
\mathfrak{P}^{(2),i,\bb{\mathfrak{P}}-1}
\right)
\right)
&=
f^{[1]@}_{s}\left(
\mathrm{sc}^{([1],2)}_{\boldsymbol{\mathcal{A}}^{(2)},s}\left(
\mathfrak{P}^{(2),i,\bb{\mathfrak{P}}-1}
\right)
\right)
\tag{1}
\\
&=
f^{[1]@}_{s}\left(
\mathrm{tg}^{([1],2)}_{\boldsymbol{\mathcal{A}}^{(2)},s}\left(
\mathfrak{P}^{(2),0,i-1}
\right)
\right)
\tag{2}
\\
&=
\mathrm{tg}^{([1],2)}_{\boldsymbol{\mathcal{B}}^{(2)},\varphi(s)}\left(
f_{s}^{(2)\flat}\left(
\mathfrak{P}^{(2),0,i-1}
\right)
\right).
\tag{3}
\end{align*}

The first equality follows from the assumption of the Artinian recursion; the second equality follows from Proposition~\ref{PDPthRecons}; finally, the last equality follows from the assumption of the Artinian recursion.

In this case, we set
\allowdisplaybreaks
\begin{align*}
f^{(2)\flat}_{s}\left(
\mathfrak{P}^{(2)}
\right)
&=
f^{(2)\flat}_{s}\left(
\mathfrak{P}^{(2),i,\bb{\mathfrak{P}^{(2)}}-1}
\right)
\circ_{\varphi(s)}^{1\mathbf{Pth}_{\boldsymbol{\mathcal{B}}^{(2)}}}
f^{(2)\flat}_{s}\left(
\mathfrak{P}^{(2),0,i-1}
\right).
\end{align*}
The above chain of equalities justifies that such composite exists.

This finishes the definition of the value of the mapping $f^{(2)\flat}$ at a path of length strictly greater than one containing at least one echelon.

If~(2), i.e., if $\mathfrak{P}^{(2)}$ is an echelonless second-order path in $\mathrm{Pth}_{\boldsymbol{\mathcal{A}}^{(2)},s}$. It could be the case that $(2.1)$ $\mathfrak{P}^{(2)}$ is not head-constant. Then let $i \in \bb{\mathfrak{P}^{(2)}}$ be the maximum index for which the subpath $\mathfrak{P}^{(2),0,i}$ of $\mathfrak{P}^{(2)}$ is a head-constant, echelonless second-order path. Note that the pairs $(\mathfrak{P}^{(2),0,i}, s)$ and $(\mathfrak{P}^{(2), i+1, \bb{\mathfrak{P}^{(2)}}-1}, s)$ $\prec_{\mathbf{Pth}_{\boldsymbol{\mathcal{A}}^{(2)}}}$-precede the pair $(\mathfrak{P}^{(2)},s)$. Therefore, the values of the mapping $f^{(2)\flat}$ at $\mathfrak{P}^{(2),0,i}$ and $\mathfrak{P}^{(2),i+1,\bb{\mathfrak{P}^{(2)}}-1}$, respectively, have already been defined. In particular, the following chain of equalities holds
\allowdisplaybreaks
\begin{align*}
\mathrm{sc}^{([1],2)}_{\boldsymbol{\mathcal{B}}^{(2)},\varphi(s)}\left(
f_{s}^{(2)\flat}\left(
\mathfrak{P}^{(2),i+1,\bb{\mathfrak{P}^{(2)}}-1}
\right)
\right)
&=
f^{[1]@}_{s}\left(
\mathrm{sc}^{([1],2)}_{\boldsymbol{\mathcal{A}}^{(2)},s}\left(
\mathfrak{P}^{(2),i+1,\bb{\mathfrak{P}^{(2)}}-1}
\right)
\right)
\tag{1}
\\
&=
f^{[1]@}_{s}\left(
\mathrm{tg}^{([1],2)}_{\boldsymbol{\mathcal{A}}^{(2)},s}\left(
\mathfrak{P}^{(2),0,i}
\right)
\right)
\tag{2}
\\
&=
\mathrm{tg}^{([1],2)}_{\boldsymbol{\mathcal{B}}^{(2)},\varphi(s)}\left(
f_{s}^{(2)\flat}\left(
\mathfrak{P}^{(2),0,i}
\right)
\right).
\tag{3}
\end{align*}

The first equality follows from the assumption of the Artinian recursion;
the second equality follows from Proposition~\ref{PDPthRecons};
finally, the last equality follows from the assumption of the Artinian recursion.

In this case, we set
\allowdisplaybreaks
\begin{align*}
f^{(2)\flat}_{s}\left(
\mathfrak{P}^{(2)}
\right)
&=
f^{(2)\flat}_{s}\left(
\mathfrak{P}^{(2),i+1,\bb{\mathfrak{P}^{(2)}}-1}
\right)
\circ_{\varphi(s)}^{1\mathbf{Pth}_{\boldsymbol{\mathcal{B}}^{(2)}}}
f^{(2)\flat}_{s}\left(
\mathfrak{P}^{(2),0,i}
\right).
\end{align*}
The above chain of equalities justifies that such composite exists.

Therefore we are left with the case of $\mathfrak{P}^{(2)}$ being a head-constant echelonless second-order path. It could be the case that $(2.2)$ $\mathfrak{P}^{(2)}$ is not coherent. Then let $i \in \bb{\mathfrak{P}^{(2)}}$ be the maximum index for which the subpath $\mathfrak{P}^{(2),0,i}$ of $\mathfrak{P}^{(2)}$ is a coherent head-constant echeloneless second-order path. Note that the pairs $(\mathfrak{P}^{(2),0,i}, s)$ and $(\mathfrak{P}^{(2), i+1, \bb{\mathfrak{P}^{(2)}-1}}, s)$ $\prec_{\mathbf{Pth}_{\boldsymbol{\mathcal{A}}^{(2)}}}$-precede the pair $(\mathfrak{P}^{(2)},s)$. Therefore, the values of the mapping $f^{(2)\flat}$ at $\mathfrak{P}^{(2),0,i}$ and $\mathfrak{P}^{(2),i+1,\bb{\mathfrak{P}}-1}$, respectively, have already been defined. In particular, the following chain of equalities holds
\allowdisplaybreaks
\begin{align*}
\mathrm{sc}^{([1],2)}_{\boldsymbol{\mathcal{B}}^{(2)},\varphi(s)}\left(
f_{s}^{(2)\flat}\left(
\mathfrak{P}^{(2),i+1,\bb{\mathfrak{P}^{(2)}}-1}
\right)
\right)
&=
f^{[1]@}_{s}\left(
\mathrm{sc}^{([1],2)}_{\boldsymbol{\mathcal{A}}^{(2)},s}\left(
\mathfrak{P}^{(2),i+1,\bb{\mathfrak{P}^{(2)}}-1}
\right)
\right)
\tag{1}
\\
&=
f^{[1]@}_{s}\left(
\mathrm{tg}^{([1],2)}_{\boldsymbol{\mathcal{A}}^{(2)},s}\left(
\mathfrak{P}^{(2),0,i}
\right)
\right)
\tag{2}
\\
&=
\mathrm{tg}^{([1],2)}_{\boldsymbol{\mathcal{B}}^{(2)},\varphi(s)}\left(
f_{s}^{(2)\flat}\left(
\mathfrak{P}^{(2),0,i}
\right)
\right).
\tag{3}
\end{align*}

The first equality follows from the assumption of the Artinian recursion;
the second equality follows from Proposition~\ref{PDPthRecons};
finally, the last equality follows from the assumption of the Artinian recursion.

In this case, we set
\allowdisplaybreaks
\begin{align*}
f^{(2)\flat}_{s}\left(
\mathfrak{P}^{(2)}
\right)
&=
f^{(2)\flat}_{s}\left(
\mathfrak{P}^{(2),i+1,\bb{\mathfrak{P}}-1}
\right)
\circ_{\varphi(s)}^{1\mathbf{Pth}_{\boldsymbol{\mathcal{B}}^{(2)}}}
f^{(2)\flat}_{s}\left(
\mathfrak{P}^{(2),0,i}
\right).
\end{align*}
The above chain of equalities justifies that such composite exists.

Therefore we are left with the case $(2.3)$ of $\mathfrak{P}^{(2)}$ being a coherent head-constant echelonless second-order path. Under this setting, the conditions for the second-order extraction algorithm, that is, Lemma~\ref{LDPthExtract}, are fulfilled. Then there exists a unique word $\mathbf{s} \in S^{\star}-\{\lambda\}$ and a unique operation symbol $\tau \in \Sigma^{\boldsymbol{\mathcal{A}}^{(1)}}_{\mathbf{s}, s}$ associated to $\mathfrak{P}^{(2)}$. Let $(\mathfrak{P}^{(2)}_{j})_{j \in \bb{\mathbf{s}}}$ be the family of second-order paths in $\mathrm{Pth}_{\boldsymbol{\mathcal{A}}^{(2)}, \mathbf{s}}$ which, in virtue of Lemma~\ref{LDPthExtract}, we can extract from $\mathfrak{P}^{(2)}$. Note that, for every $j \in \bb{\mathbf{s}}$, we have that $(\mathfrak{P}^{(2)}_{j}, s_{j}) \prec_{\mathbf{Pth}_{\boldsymbol{\mathcal{A}}^{(2)}}} (\mathfrak{P}^{(2)}, s)$. Therefore, for every $j \in \bb{\mathbf{s}}$, the value of the mapping $f^{(2)\flat}$ at $\mathfrak{P}_{j}$ has already been defined.

It could be the case that $(2.3.1)$ $\tau$ is $\sigma$ an operation symbol in $\Sigma_{\mathbf{s}, s}$, that is, $\mathfrak{P}^{(2)}$ is $\sigma^{\mathbf{Pth}_{\boldsymbol{\mathcal{A}}^{(2)}}}((\mathfrak{P}^{(2)}_{j})_{j\in\bb{\mathbf{s}}})$.

In this case, we set
$$
f^{(2)\flat}_{s}\left(
\mathfrak{P}^{(2)}
\right)
=
\sigma^{\mathbf{Pth}_{\boldsymbol{\mathcal{B}}^{(2)}}^{\mathbf{f}^{(2)}(0,2)}}\left(\left(
f^{(2)\flat}_{s_{j}}\left(
\mathfrak{P}^{(2)}_{j}
\right)
\right)_{j\in\bb{\mathbf{s}}}\right).
$$

Finally, we are left with case $(2.3.2)$ $\tau$ is the operation symbol $\circ_{s}^{0}$ in $\Sigma^{\boldsymbol{\mathcal{A}}^{(1)}}_{ss,s}$, that is, $\mathfrak{P}^{(2)}$ is $\mathfrak{P}^{(2)}_{1} \circ_{s}^{0\mathbf{Pth}_{\boldsymbol{\mathcal{A}}^{(2)}}} \mathfrak{P}^{(2)}_{0}$. Note that, since the $0$-composition is well-defined,
$$
\mathrm{tg}_{\boldsymbol{\mathcal{A}}^{(2)}, s}^{(0,2)} \left(
\mathfrak{P}^{(2)}_{0}
\right)
=
\mathrm{sc}_{\boldsymbol{\mathcal{A}}^{(2)}, s}^{(0,2)} \left(
\mathfrak{P}^{(2)}_{1}
\right).
$$

Let us also note that the following chain of equalities holds
\begin{flushleft}
$
\mathrm{tg}_{\boldsymbol{\mathcal{B}}^{(2)}, \varphi(s)}^{(0,2)} \left(
f^{(2)\flat}_{s}\left(
\mathfrak{P}^{(2)}_{0}
\right)
\right)
$
\allowdisplaybreaks
\begin{align*}
&=
\mathrm{tg}_{\boldsymbol{\mathcal{B}}^{(2)}, \varphi(s)}^{(0,[1])} \left(
\mathrm{ip}^{([1],Y)@}_{\boldsymbol{\mathcal{B}}^{(2)}, \varphi(s)}\left(
\mathrm{tg}_{\boldsymbol{\mathcal{B}}^{(2)}, \varphi(s)}^{([1],2)} \left(
f^{(2)\flat}_{s}\left(
\mathfrak{P}^{(2)}_{0}
\right)
\right)
\right)
\right)
\tag{1}
\\
&=
\mathrm{tg}_{\boldsymbol{\mathcal{B}}^{(2)}, \varphi(s)}^{(0,[1])} \left(
\mathrm{ip}^{([1],Y)@}_{\boldsymbol{\mathcal{B}}^{(2)}, \varphi(s)}\left(
f^{[1]@}_{s}\left(
\mathrm{tg}_{\boldsymbol{\mathcal{A}}^{(2)}, s}^{([1],2)} \left(
\mathfrak{P}^{(2)}_{0}
\right)
\right)
\right)
\right)
\tag{2}
\\
&=
\mathrm{tg}_{\boldsymbol{\mathcal{B}}^{(2)}, \varphi(s)}^{(0,[1])} \left(
f^{[1]@}_{s}\left(
\mathrm{ip}^{([1],X)@}_{\boldsymbol{\mathcal{A}}^{(2)}, s}\left(
\mathrm{tg}_{\boldsymbol{\mathcal{A}}^{(2)}, s}^{([1],2)} \left(
\mathfrak{P}^{(2)}_{0}
\right)
\right)
\right)
\right)
\tag{3}
\\
&=
f^{(0)\sharp}_{s}\left(
\mathrm{tg}_{\boldsymbol{\mathcal{A}}^{(2)}, s}^{(0,[1])} \left(
\mathrm{ip}^{([1],X)@}_{\boldsymbol{\mathcal{A}}^{(2)}, s}\left(
\mathrm{tg}_{\boldsymbol{\mathcal{A}}^{(2)}, s}^{([1],2)} \left(
\mathfrak{P}^{(2)}_{0}
\right)
\right)
\right)
\right)
\tag{4}
\\
&=
f^{(0)\sharp}_{s}\left(
\mathrm{tg}^{(0,2)}_{\boldsymbol{\mathcal{A}}^{(2)}, s}\left(
\mathfrak{P}^{(2)}_{0}
\right)
\right)
\tag{5}
\\
&=
f^{(0)\sharp}_{s}\left(
\mathrm{sc}^{(0,2)}_{\boldsymbol{\mathcal{A}}^{(2)}, s}\left(
\mathfrak{P}^{(2)}_{1}
\right)
\right)
\tag{6}
\\
&=
f^{(0)\sharp}_{s}\left(
\mathrm{sc}_{\boldsymbol{\mathcal{A}}^{(2)}, s}^{(0,[1])} \left(
\mathrm{ip}^{([1],X)@}_{\boldsymbol{\mathcal{A}}^{(2)}, s}\left(
\mathrm{sc}_{\boldsymbol{\mathcal{A}}^{(2)}, s}^{([1],2)} \left(
\mathfrak{P}^{(2)}_{1}
\right)
\right)
\right)
\right)
\tag{7}
\\
&=
\mathrm{sc}_{\boldsymbol{\mathcal{B}}^{(2)}, \varphi(s)}^{(0,[1])} \left(
f^{[1]@}_{s}\left(
\mathrm{ip}^{([1],X)@}_{\boldsymbol{\mathcal{A}}^{(2)}, s}\left(
\mathrm{sc}_{\boldsymbol{\mathcal{A}}^{(2)}, s}^{([1],2)} \left(
\mathfrak{P}^{(2)}_{1}
\right)
\right)
\right)
\right)
\tag{8}
\\
&=
\mathrm{sc}_{\boldsymbol{\mathcal{B}}^{(2)}, \varphi(s)}^{(0,[1])} \left(
\mathrm{ip}^{([1],Y)@}_{\boldsymbol{\mathcal{B}}^{(2)}, \varphi(s)}\left(
f^{[1]@}_{s}\left(
\mathrm{sc}_{\boldsymbol{\mathcal{A}}^{(2)}, s}^{([1],2)} \left(
\mathfrak{P}^{(2)}_{1}
\right)
\right)
\right)
\right)
\tag{9}
\\
&=
\mathrm{sc}_{\boldsymbol{\mathcal{B}}^{(2)}, \varphi(s)}^{(0,[1])} \left(
\mathrm{ip}^{([1],Y)@}_{\boldsymbol{\mathcal{B}}^{(2)}, \varphi(s)}\left(
\mathrm{sc}_{\boldsymbol{\mathcal{B}}^{(2)}, \varphi(s)}^{([1],2)} \left(
f^{(2)\flat}_{s}\left(
\mathfrak{P}^{(2)}_{1}
\right)
\right)
\right)
\right)
\tag{10}
\\
&=
\mathrm{sc}_{\boldsymbol{\mathcal{B}}^{(2)}, \varphi(s)}^{(0,2)} \left(
f^{(2)\flat}_{s}\left(
\mathfrak{P}^{(2)}_{1}
\right)
\right).
\tag{11}
\end{align*}
\end{flushleft}

The first equality unravels the definition of the mapping $\mathrm{tg}^{(0,2)}_{\boldsymbol{\mathcal{B}}^{(2)}}$, introduced in Definition~\ref{DDScTgZ};
the second equality follows from the assumption of the Artinian recursion;
the third equality follows from Proposition~\ref{PQPthExtIp};
the fourth equality follows from Proposition~\ref{PQPthExtScTg};
the fifth equality recovers the definition of the mapping $\mathrm{tg}^{(0,2)}_{\boldsymbol{\mathcal{A}}^{(2)}}$, introduced in Definition~\ref{DDScTgZ};
the sixth equality follows from the fact that, since the $0$-composition is well-defined,
$
\mathrm{tg}_{\boldsymbol{\mathcal{A}}^{(2)}, s}^{(0,2)} (\mathfrak{P}^{(2)}_{0})
=
\mathrm{sc}_{\boldsymbol{\mathcal{A}}^{(2)}, s}^{(0,2)} (\mathfrak{P}^{(2)}_{1});
$
the seventh equality unravels the definition of the mapping $\mathrm{sc}^{(0,2)}_{\boldsymbol{\mathcal{A}}^{(2)}}$, introduced in Definition~\ref{DDScTgZ};
the eighth equality follows from Proposition~\ref{PPthExt};
the ninth equality follows from Proposition~\ref{PQPthExtIp};
the tenth equality follows from the assumption of the Artinian recursion;
finally, the last equality recovers the definition of the mapping $\mathrm{sc}^{(0,2)}_{\boldsymbol{\mathcal{B}}^{(2)}}$, introduced in Definition~\ref{DDScTgZ};

In this case we set
$$
f ^{(2)\flat}_{s}\left(
\mathfrak{P}^{(2)}
\right)
=
f^{(2)\flat}_{s}\left(
\mathfrak{P}^{(2)}_{1}
\right)
\circ_{\varphi(s)}^{0\mathbf{Pth}_{\boldsymbol{\mathcal{B}}^{(2)}}}
f^{(2)\flat}_{s}\left(
\mathfrak{P}^{(2)}_{0}
\right)
$$
The above chain of equalities justifies that such composite exists.

That finishes the definition of the value of the mapping $f^{(2)\flat}$ at an echelonless second-order path.

This completes the definition of the mapping $f^{(2)\flat}$.

Now, we prove that the just defined mapping satisfies all the properties listed in the proposition.

\textsf{(1)}
$\mathrm{sc}^{([1],2)}_{\boldsymbol{\mathcal{B}}^{(2)},\varphi}\circ f^{(2)\flat}=f^{[1]@}\circ \mathrm{sc}^{([1],2)}_{\boldsymbol{\mathcal{A}}^{(2)}}$.

Let $s$ be a sort in $S$ and let $\mathfrak{P}^{(2)}$ be a second-order path in $\mathrm{Pth}_{\boldsymbol{\mathcal{A}}^{(2)}, s}$. We prove that
$$
\mathrm{sc}^{([1],2)}_{\boldsymbol{\mathcal{B}}^{(2)},\varphi(s)}\left(
f_{s}^{(2)\flat} \left(
\mathfrak{P}^{(2)}
\right)
\right)
=
f_{s}^{[1]@}\left(
\mathrm{sc}^{([1],2)}_{\boldsymbol{\mathcal{A}}^{(2)}, s}\left(
\mathfrak{P}^{(2)}
\right)
\right)
$$
by Artinian induction on $(\coprod \mathrm{Pth}_{\boldsymbol{\mathcal{A}}^{(2)}}, \leq_{\mathbf{Pth}_{\boldsymbol{\mathcal{A}}^{(2)}}})$.

{\sffamily Base step of the Artinian induction.}

Let $(\mathfrak{P}^{(2)}, s)$ be a minimal element of $(\coprod \mathrm{Pth}_{\boldsymbol{\mathcal{A}}^{(2)}}, \leq_{\mathbf{Pth}_{\boldsymbol{\mathcal{A}}^{(2)}}})$. Then, by Proposition~\ref{PDMinimal}, the path $\mathfrak{P}^{(2)}$ is either~(1) an $(2,[1])$-identity second-order path or~(2) a second-order echelon.

If~(1), i.e., if $\mathfrak{P}^{(2)}$ is an $(2,[1])$-identity second-order path, then $\mathfrak{P}^{(2)}=\mathrm{ip}^{(2,[1])\sharp}_{\boldsymbol{\mathcal{A}}^{(2)}, s}([P]_{s})$ for some term $[P]_{s} \in [\mathrm{PT}_{\boldsymbol{\mathcal{A}}^{(2)}}]_{s}$. Note that the following chain of equalities holds
\allowdisplaybreaks
\begin{align*}
\mathrm{sc}_{\boldsymbol{\mathcal{B}}^{(2)},\varphi(s)}^{([1],2)}\left(
f^{(2)\flat}_{s}\left(
\mathfrak{P}^{(2)}
\right)
\right) 
&= 
\mathrm{sc}_{\boldsymbol{\mathcal{B}}^{(2)},\varphi(s)}^{([1],2)}\left(
\mathrm{ip}_{\boldsymbol{\mathcal{B}}^{(2)},\varphi(s)}^{(2,[1])\sharp}\left(
f_{s}^{[1]@}\left(
\left[P\right]_{s}
\right)
\right)
\right)
\tag{1}
\\
&=
f_{s}^{[1]@}\left(
\left[P\right]_{s}
\right)
\tag{2}
\\
&=
f_{s}^{[1]@}\left(
\mathrm{sc}_{\boldsymbol{\mathcal{A}}^{(2)},s}^{([1],2)}\left(
\mathrm{ip}_{\boldsymbol{\mathcal{A}}^{(2)},s}^{(2,[1])\sharp}\left(
\left[P\right]_{s}
\right)
\right)
\right)
\tag{3}
\\
&= f^{[1]@}_{s}\left(
\mathrm{sc}_{\boldsymbol{\mathcal{A}}^{(2)},s}^{([1],2)}
\left(
\mathfrak{P}^{(2)}
\right)
\right).
\tag{4}
\end{align*}

The first equality unravels the definition of $f_{s}^{(2)\flat}\left(\mathfrak{P}^{(2)}\right)$;
the second and third equalities follows from the fact that, according to Proposition~\ref{PDBasicEq},
\begin{align*}
\mathrm{sc}_{\boldsymbol{\mathcal{B}}^{(1)},\varphi(s)}^{([1],2)} \circ \mathrm{ip}^{(2,[1])\sharp}_{\boldsymbol{\mathcal{B}}^{(2)}, \varphi(s)} 
&=
\mathrm{id}^{[\mathrm{PT}_{\boldsymbol{\mathcal{B}}^{(1)}}]_{\varphi(s)}}
&&\mbox{and}&
\mathrm{sc}_{\boldsymbol{\mathcal{A}}^{(1)},s}^{(0,1)} \circ \mathrm{ip}^{(1,0)\sharp}_{\boldsymbol{\mathcal{A}}^{(1)}, s} 
&= 
\mathrm{id}^{[\mathrm{PT}_{\boldsymbol{\mathcal{A}}^{(1)}}]_{s}};
\end{align*}
the fourth equality recovers the definition of $\mathfrak{P}^{(2)}$.

If~(2), i.e., if $\mathfrak{P}^{(2)}$ is an echelon associated to a rewrite rule $\mathfrak{p}^{(2)}=([M]_{s}, [N]_{s})$, the following chain of equalities holds
\allowdisplaybreaks
\begin{align*}
\mathrm{sc}_{\boldsymbol{\mathcal{B}}^{(2)},\varphi(s)}^{([1],2)}\left(
f^{(2)\flat}_{s}\left(
\mathfrak{P}
\right)
\right)
&=
\mathrm{sc}_{\boldsymbol{\mathcal{B}}^{(2)},\varphi(s)}^{([1],2)}\left(
f^{(1)}_{s}\left(
\mathfrak{p}^{(2)}
\right)
\right) 
\tag{1}
\\
&=
f^{[1]@}\left(
[M]_{s}
\right)
\tag{2}
\\
&=
f^{[1]@}\left(
\mathrm{sc}_{\boldsymbol{\mathcal{A}}^{(2)}, s}^{([1],2)} \left(
\mathfrak{P}^{(2)}
\right)
\right).
\tag{3}
\end{align*}

The first equality unravels the definition of $f_{s}^{(2)\flat}\left(\mathfrak{P}^{(2)}\right)$;
the second equality follows from item~(2) of Definition~\ref{DRws2Mor};
finally, the last equality follows from Proposition~\ref{PDBasicEqA}.

This completes the base step of the Artinian induction.

{\sffamily Inductive step of the Artinian induction.}

Let $(\mathfrak{P}^{(2)},s)$ be a non-minimal element of $(\coprod\mathrm{Pth}_{\boldsymbol{\mathcal{A}}^{(2)}}, \leq_{\mathbf{Pth}_{\boldsymbol{\mathcal{A}}^{(2)}}})$. We can assume that $\mathfrak{P}^{(2)}$ is a not a $(2,[1])$-identity second-order path, since case has already been proven. Let us suppose that, for every sort $t\in S$ and every second-order path $\mathfrak{Q}^{(2)}\in\mathrm{Pth}_{\boldsymbol{\mathcal{A}}^{(1)},t}$, if $(\mathfrak{Q}^{(2)},t)<_{\mathbf{Pth}_{\boldsymbol{\mathcal{A}}^{(1)}}}(\mathfrak{P}^{(2)},s)$, then the following equality holds
$$
\mathrm{sc}_{\boldsymbol{\mathcal{B}}^{(2)},\varphi(t)}^{([1],2)}\left(f^{(2)\flat}_{t}\left(\mathfrak{Q}^{(2)}\right)\right) 
=
f^{[1]@}_{t}\left(\mathrm{sc}_{\boldsymbol{\mathcal{A}}^{(2)},t}^{([1],2)}\left(\mathfrak{Q}^{(2)}\right)\right)
$$

By Lemma~\ref{LDOrdI}, we have that $\mathfrak{P}^{(2)}$ is either~(1) a path of length strictly greater than one containing at least one second-order echelon or~(2) an echelonless second-order path.

If~(1), i.e., if $\mathfrak{P}^{(2)}$ is a path of length strictly greater than one containing at least one second-order echelon, then let $i\in \bb{\mathfrak{P}^{(2)}}$ be the first index for which the one-step subpath $\mathfrak{P}^{(2),i,i}$ of $\mathfrak{P}^{(2)}$ is a second-order echelon. We consider different cases for $i$ according to the cases presented in Definition~\ref{DDOrd}.

If $i=0$, we have that the pairs $(\mathfrak{P}^{(2),0,0},s)$ and $(\mathfrak{P}^{(2),1,\bb{\mathfrak{P}^{(2)}}-1},s)$ $\prec_{\mathbf{Pth}_{\boldsymbol{\mathcal{A}}^{(2)}}}$-precede the pair $(\mathfrak{P}^{(2)},s)$. The following chain of equalities holds
\allowdisplaybreaks
\begin{flushleft}
$
\mathrm{sc}_{\boldsymbol{\mathcal{B}}^{(2)},\varphi(s)}^{([1],2)}\left(
f_{s}^{(2)\flat}\left(
\mathfrak{P}^{(2)}
\right)
\right)
$
\begin{align*}
&=
\mathrm{sc}_{\boldsymbol{\mathcal{B}}^{(2)},\varphi(s)}^{([1],2)}\left(
f^{(2)\flat}_{s}\left(
\mathfrak{P}^{(2),1,\bb{\mathfrak{P}}-1}
\right)
\circ_{\varphi(s)}^{1\mathbf{Pth}_{\boldsymbol{\mathcal{B}}^{(2)}}}
f^{(2)\flat}_{s}\left(
\mathfrak{P}^{(2),0,0}
\right)
\right)
\tag{1}
\\
&=
\mathrm{sc}_{\boldsymbol{\mathcal{B}}^{(2)},\varphi(s)}^{([1],2)}\left(
f^{(2)\flat}_{s}\left(
\mathfrak{P}^{(2),0,0}
\right)
\right)
\tag{2}
\\
&=
f^{[1]@}_{s}\left(
\mathrm{sc}_{\boldsymbol{\mathcal{A}}^{(2)},s}^{([1],2)}\left(
\mathfrak{P}^{(2),0,0}
\right)
\right)
\tag{3}
\\
&=
f^{[1]@}_{s}\left(
\mathrm{sc}_{\boldsymbol{\mathcal{A}}^{(2)},s}^{([1],2)}\left(
\mathfrak{P}^{(2),1,\bb{\mathfrak{P}}-1}
\circ_{s}^{1\mathbf{Pth}_{\boldsymbol{\mathcal{A}}^{(2)}}}
\mathfrak{P}^{(2),0,0}
\right)
\right)
\tag{4}
\\
&=
f^{[1]@}_{s}\left(
\mathrm{sc}_{\boldsymbol{\mathcal{A}}^{(2)},s}^{([1],2)}\left(
\mathfrak{P}^{(2)}
\right)
\right).
\tag{5}
\end{align*}
\end{flushleft}

The first equality unravels the definition of $f_{s}^{(2)\flat}\left(\mathfrak{P}^{(2)}\right)$;
the second equality follows from Proposition~\ref{PDPthComp};
the third equality follows by Artinian induction;
the fourth equality follows from Proposition~\ref{PDPthComp};
the fifth equality recovers the definition of the path $\mathfrak{P}^{(2)}$.

If $i\neq 0$, we have that the pairs $(\mathfrak{P}^{(2),0,i-1},s)$ and $(\mathfrak{P}^{(2),i,\bb{\mathfrak{P}^{(2)}}-1},s)$ $\prec_{\mathbf{Pth}_{\boldsymbol{\mathcal{A}}^{(2)}}}$-precede the pair $(\mathfrak{P}^{(2)},s)$. The following chain of equalities holds
\allowdisplaybreaks
\begin{flushleft}
$
\mathrm{sc}_{\boldsymbol{\mathcal{B}}^{(2)},\varphi(s)}^{([1],2)}\left(
f_{s}^{(2)\flat}\left(
\mathfrak{P}^{(2)}
\right)
\right)
$
\begin{align*}
&=
\mathrm{sc}_{\boldsymbol{\mathcal{B}}^{(2)},\varphi(s)}^{([1],2)}\left(
f^{(2)\flat}_{s}\left(
\mathfrak{P}^{(2),i,\bb{\mathfrak{P}^{(2)}}-1}
\right)
\circ_{\varphi(s)}^{1\mathbf{Pth}_{\boldsymbol{\mathcal{B}}^{(2)}}}
f^{(2)\flat}_{s}\left(
\mathfrak{P}^{(2),0,i-1}
\right)
\right)
\tag{1}
\\
&=
\mathrm{sc}_{\boldsymbol{\mathcal{B}}^{(2)},\varphi(s)}^{([1],2)}\left(
f^{(2)\flat}_{s}\left(
\mathfrak{P}^{(2),0,i-1}
\right)
\right)
\tag{2}
\\
&=
f^{[1]@}_{s}\left(
\mathrm{sc}_{\boldsymbol{\mathcal{A}}^{(2)},s}^{([1],2)}\left(
\mathfrak{P}^{(2),0,1-i}
\right)
\right)
\tag{3}
\\
&=
f^{[1]@}_{s}\left(
\mathrm{sc}_{\boldsymbol{\mathcal{A}}^{(2)},s}^{([1],2)}\left(
\mathfrak{P}^{(2),i,\bb{\mathfrak{P}^{(2)}}-1}
\circ_{s}^{1\mathbf{Pth}_{\boldsymbol{\mathcal{A}}^{(2)}}}
\mathfrak{P}^{(2),0,i-1}
\right)
\right)
\tag{4}
\\
&=
f^{[1]@}_{s}\left(
\mathrm{sc}_{\boldsymbol{\mathcal{A}}^{(2)},s}^{([1],2)}\left(
\mathfrak{P}^{(2)}
\right)
\right).
\tag{5}
\end{align*}
\end{flushleft}

The first equality unravels the definition of $f_{s}^{(2)\flat}\left(\mathfrak{P}^{(2)}\right)$;
the second equality follows from Proposition~\ref{PDPthComp};
the third equality follows by Artinian induction;
the fourth equality follows from Proposition~\ref{PDPthComp};
the fifth equality recovers the definition of the path $\mathfrak{P}^{(2)}$.

Similarly, the following equality holds
$$
\mathrm{tg}_{\boldsymbol{\mathcal{B}}^{(2)},\varphi(s)}^{([1],2)}\left(
f_{s}^{(2)\flat}\left(
\mathfrak{P}^{(2)}
\right)
\right)
=
f^{[1]@}_{s}\left(
\mathrm{tg}_{\boldsymbol{\mathcal{A}}^{(2)},s}^{([1],2)}\left(
\mathfrak{P}^{(2)}
\right)
\right).
$$

If~(2), i.e., if $\mathfrak{P}^{(2)}$ is an echelonless second-order path in $\mathrm{Pth}_{\boldsymbol{\mathcal{A}}^{(2)},s}$. It could be the case that $(2.1)$ $\mathfrak{P}^{(2)}$ is not head-constant. Then let $i \in \bb{\mathfrak{P}^{(2)}}$ be the maximum index for which the subpath $\mathfrak{P}^{(2),0,i}$ of $\mathfrak{P}^{(2)}$ is a head-constant, echelonless second-order path. Note that the pairs $(\mathfrak{P}^{(2),0,i}, s)$ and $(\mathfrak{P}^{(2), i+1, \bb{\mathfrak{P}^{(2)}}-1}, s)$ $\prec_{\mathbf{Pth}_{\boldsymbol{\mathcal{A}}^{(2)}}}$-precede the pair $(\mathfrak{P}^{(2)},s)$. Thus, the following chain of equalities holds
\allowdisplaybreaks
\begin{flushleft}
$
\mathrm{sc}_{\boldsymbol{\mathcal{B}}^{(2)},\varphi(s)}^{([1],2)}\left(
f_{s}^{(2)\flat}\left(
\mathfrak{P}^{(2)}
\right)
\right)
$
\begin{align*}
&=
\mathrm{sc}_{\boldsymbol{\mathcal{B}}^{(2)},\varphi(s)}^{([1],2)}\left(
f^{(2)\flat}_{s}\left(
\mathfrak{P}^{(2),i+1,\bb{\mathfrak{P}^{(2)}}-1}
\right)
\circ_{\varphi(s)}^{1\mathbf{Pth}_{\boldsymbol{\mathcal{B}}^{(2)}}}
f^{(2)\flat}_{s}\left(
\mathfrak{P}^{(2),0,i}
\right)
\right)
\tag{1}
\\
&=
\mathrm{sc}_{\boldsymbol{\mathcal{B}}^{(2)},\varphi(s)}^{([1],2)}\left(
f^{(2)\flat}_{s}\left(
\mathfrak{P}^{(2),0,i}
\right)
\right)
\tag{2}
\\
&=
f^{[1]@}_{s}\left(
\mathrm{sc}_{\boldsymbol{\mathcal{A}}^{(2)},s}^{([1],2)}\left(
\mathfrak{P}^{(2),0,1}
\right)
\right)
\tag{3}
\\
&=
f^{[1]@}_{s}\left(
\mathrm{sc}_{\boldsymbol{\mathcal{A}}^{(2)},s}^{([1],2)}\left(
\mathfrak{P}^{(2),i+1,\bb{\mathfrak{P}^{(2)}}-1}
\circ_{s}^{1\mathbf{Pth}_{\boldsymbol{\mathcal{A}}^{(2)}}}
\mathfrak{P}^{(2),0,i}
\right)
\right)
\tag{4}
\\
&=
f^{[1]@}_{s}\left(
\mathrm{sc}_{\boldsymbol{\mathcal{A}}^{(2)},s}^{([1],2)}\left(
\mathfrak{P}^{(2)}
\right)
\right).
\tag{5}
\end{align*}
\end{flushleft}

The first equality unravels the definition of $f_{s}^{(2)\flat}\left(\mathfrak{P}^{(2)}\right)$;
the second equality follows from Proposition~\ref{PDPthComp};
the third equality follows by Artinian induction;
the fourth equality follows from Proposition~\ref{PDPthComp};
the fifth equality recovers the definition of the path $\mathfrak{P}^{(2)}$.

Therefore we are left with the case of $\mathfrak{P}^{(2)}$ being a head-constant echelonless second-order path. It could be the case that $(2.2)$ $\mathfrak{P}^{(2)}$ is not coherent. Then let $i \in \bb{\mathfrak{P}^{(2)}}$ be the maximum index for which the subpath $\mathfrak{P}^{(2),0,i}$ of $\mathfrak{P}^{(2)}$ is a coherent head-constant echeloneless second-order path. Note that the pairs $(\mathfrak{P}^{(2),0,i}, s)$ and $(\mathfrak{P}^{(2), i+1, \bb{\mathfrak{P}^{(2)}-1}}, s)$ $\prec_{\mathbf{Pth}_{\boldsymbol{\mathcal{A}}^{(2)}}}$-precede the pair $(\mathfrak{P}^{(2)},s)$. Thus, the following chain of equalities holds
\allowdisplaybreaks
\begin{flushleft}
$
\mathrm{sc}_{\boldsymbol{\mathcal{B}}^{(2)},\varphi(s)}^{([1],2)}\left(
f_{s}^{(2)\flat}\left(
\mathfrak{P}^{(2)}
\right)
\right)
$
\begin{align*}
&=
\mathrm{sc}_{\boldsymbol{\mathcal{B}}^{(2)},\varphi(s)}^{([1],2)}\left(
f^{(2)\flat}_{s}\left(
\mathfrak{P}^{(2),i+1,\bb{\mathfrak{P}^{(2)}}-1}
\right)
\circ_{\varphi(s)}^{1\mathbf{Pth}_{\boldsymbol{\mathcal{B}}^{(2)}}}
f^{(2)\flat}_{s}\left(
\mathfrak{P}^{(2),0,i}
\right)
\right)
\tag{1}
\\
&=
\mathrm{sc}_{\boldsymbol{\mathcal{B}}^{(2)},\varphi(s)}^{([1],2)}\left(
f^{(2)\flat}_{s}\left(
\mathfrak{P}^{(2),0,i}
\right)
\right)
\tag{2}
\\
&=
f^{[1]@}_{s}\left(
\mathrm{sc}_{\boldsymbol{\mathcal{A}}^{(2)},s}^{([1],2)}\left(
\mathfrak{P}^{(2),0,1}
\right)
\right)
\tag{3}
\\
&=
f^{[1]@}_{s}\left(
\mathrm{sc}_{\boldsymbol{\mathcal{A}}^{(2)},s}^{([1],2)}\left(
\mathfrak{P}^{(2),i+1,\bb{\mathfrak{P}^{(2)}}-1}
\circ_{s}^{1\mathbf{Pth}_{\boldsymbol{\mathcal{A}}^{(2)}}}
\mathfrak{P}^{(2),0,i}
\right)
\right)
\tag{4}
\\
&=
f^{[1]@}_{s}\left(
\mathrm{sc}_{\boldsymbol{\mathcal{A}}^{(2)},s}^{([1],2)}\left(
\mathfrak{P}^{(2)}
\right)
\right).
\tag{5}
\end{align*}
\end{flushleft}

The first equality unravels the definition of $f_{s}^{(2)\flat}\left(\mathfrak{P}^{(2)}\right)$;
the second equality follows from Proposition~\ref{PDPthComp};
the third equality follows by Artinian induction;
the fourth equality follows from Proposition~\ref{PDPthComp};
the fifth equality recovers the definition of the path $\mathfrak{P}^{(2)}$.

Similarly, the following equality holds
$$
\mathrm{tg}_{\boldsymbol{\mathcal{B}}^{(2)},\varphi(s)}^{([1],2)}\left(
f_{s}^{(2)\flat}\left(
\mathfrak{P}^{(2)}
\right)
\right)
=
f^{[1]@}_{s}\left(
\mathrm{tg}_{\boldsymbol{\mathcal{A}}^{(2)},s}^{([1],2)}\left(
\mathfrak{P}^{(2)}
\right)
\right).
$$

Therefore we are left with the case $(2.3)$ of $\mathfrak{P}^{(2)}$ being a coherent head-constant echelonless second-order path. Under this setting, the conditions for the second-order extraction algorithm, that is, Lemma~\ref{LDPthExtract}, are fulfilled. Then there exists a unique word $\mathbf{s} \in S^{\star}-\{\lambda\}$ and a unique operation symbol $\tau \in \Sigma^{\boldsymbol{\mathcal{A}}^{(1)}}_{\mathbf{s}, s}$ associated to $\mathfrak{P}^{(2)}$. Let $(\mathfrak{P}^{(2)}_{j})_{j \in \bb{\mathbf{s}}}$ be the family of second-order paths in $\mathrm{Pth}_{\boldsymbol{\mathcal{A}}^{(2)}, \mathbf{s}}$ which, in virtue of Lemma~\ref{LDPthExtract}, we can extract from $\mathfrak{P}^{(2)}$. Note that, for every $j \in \bb{\mathbf{s}}$, we have that $(\mathfrak{P}^{(2)}_{j}, s_{j}) \prec_{\mathbf{Pth}_{\boldsymbol{\mathcal{A}}^{(2)}}} (\mathfrak{P}^{(2)}, s)$. 

It could be the case that $(2.3.1)$ $\tau$ is $\sigma$ an operation symbol in $\Sigma_{\mathbf{s}, s}$, that is, $\mathfrak{P}^{(2)}$ is $\sigma^{\mathbf{Pth}_{\boldsymbol{\mathcal{A}}^{(2)}}}((\mathfrak{P}^{(2)}_{j})_{j\in\bb{\mathbf{s}}})$. Then, the following chain of equalities holds
\allowdisplaybreaks
\begin{align*}
\mathrm{sc}_{\boldsymbol{\mathcal{B}}^{(2)},\varphi(s)}^{([1],2)}\left(
f_{s}^{(2)\flat}\left(
\mathfrak{P}^{(2)}
\right)
\right)
&=
\mathrm{sc}_{\boldsymbol{\mathcal{B}}^{(2)},\varphi(s)}^{([1],2)}\left(
\sigma^{\mathbf{Pth}_{\boldsymbol{\mathcal{B}}^{(2)}}^{\mathbf{f}^{(2)}(0,2)}}\left(\left(
f_{s_{j}}^{(2)\flat}\left(
\mathfrak{P}_{j}^{(2)}
\right)
\right)_{j\in\bb{\mathbf{s}}}\right)
\right)
\tag{1}
\\
&=
\sigma^{[\mathbf{PT}_{\boldsymbol{\mathcal{B}}^{(1)}}^{\mathbf{f}^{(1)}}]}\left(\left(
\mathrm{sc}_{\boldsymbol{\mathcal{B}}^{(2)},\varphi(s_{j})}^{([1],2)}\left(
f_{s_{j}}^{(2)\flat}\left(
\mathfrak{P}_{j}^{(2)}
\right)
\right)
\right)_{j\in\bb{\mathbf{s}}}\right)
\tag{2}
\\
&=
\sigma^{[\mathbf{PT}_{\boldsymbol{\mathcal{B}}^{(1)}}^{\mathbf{f}^{(1)}}]}\left(\left(
f_{s_{j}}^{[1]@}\left(
\mathrm{sc}_{\boldsymbol{\mathcal{A}}^{(2)},s_{j}}^{([1],2)}\left(
\mathfrak{P}_{j}^{(2)}
\right)
\right)
\right)_{j\in\bb{\mathbf{s}}}\right)
\tag{3}
\\
&=
f_{s}^{[1]@}\left(
\sigma^{[\mathbf{PT}_{\boldsymbol{\mathcal{A}}^{(1)}}]}\left(\left(
\mathrm{sc}_{\boldsymbol{\mathcal{A}}^{(2)},s_{j}}^{([1],2)}\left(
\mathfrak{P}_{j}^{(2)}
\right)
\right)_{j\in\bb{\mathbf{s}}}\right)
\right)
\tag{4}
\\
&=
f_{s}^{[1]@}\left(
\mathrm{sc}_{\boldsymbol{\mathcal{A}}^{(2)},s}^{([1],2)}\left(
\sigma^{\mathbf{Pth}_{\boldsymbol{\mathcal{A}}^{(2)}}}\left(\left(
\mathfrak{P}_{j}^{(2)}
\right)_{j\in\bb{\mathbf{s}}}\right)
\right)
\right)
\tag{5}
\\
&=
f_{s}^{[1]@}\left(
\mathrm{sc}_{\boldsymbol{\mathcal{A}}^{(2)},s}^{([1],2)}\left(
\mathfrak{P}^{(2)}
\right)
\right).
\tag{6}
\end{align*}

The first equality unravels the definition of $f_{s}^{(2)\flat}\left(\mathfrak{P}^{(2)}\right)$;
the second equality follows from the fact that, by Proposition~\ref{PScTgUDBCatHom}, $\mathrm{sc}_{\boldsymbol{\mathcal{B}}^{(2)}, \varphi}^{([1],2)}$ is a $\Sigma$-homomorphism from $\mathbf{Pth}^{\mathbf{f}^{(2)}(0,2)}_{\boldsymbol{\mathcal{B}}^{(1)}}$ to $[\mathbf{PT}_{\boldsymbol{\mathcal{B}}^{(1)}}^{\mathbf{f}^{(1)}(0,2)}]$;
the third equality follows by Artinian induction;
the fourth equality follows from the fact that, following Definition~\ref{DQPthExt}, $f^{[1]@}$ is a $\Sigma^{\boldsymbol{\mathcal{A}}^{(1)}}$-homomorphism from $[\mathbf{PT}_{\boldsymbol{\mathcal{A}}^{(1)}}]$ to $[\mathbf{PT}^{\mathbf{f}^{(1)}}_{\boldsymbol{\mathcal{B}}^{(1)}}]$;
the fifth equality follows from Proposition~\ref{PDPthAlg};
finally, the last equality recovers the definition of $\mathfrak{P}^{(2)}$.

Finally, we are left with case $(2.3.2)$ $\tau$ is the operation symbol $\circ_{s}^{0}$ in $\Sigma^{\boldsymbol{\mathcal{A}}^{(1)}}_{ss,s}$, that is, $\mathfrak{P}^{(2)}$ is $\mathfrak{P}^{(2)}_{1} \circ_{s}^{0\mathbf{Pth}_{\boldsymbol{\mathcal{A}}^{(2)}}} \mathfrak{P}^{(2)}_{0}$. Note that, since the $0$-composition is well-defined,
$$
\mathrm{tg}_{\boldsymbol{\mathcal{A}}^{(2)}, s}^{(0,2)} \left(
\mathfrak{P}^{(2)}_{0}
\right)
=
\mathrm{sc}_{\boldsymbol{\mathcal{A}}^{(2)}, s}^{(0,2)} \left(
\mathfrak{P}^{(2)}_{1}
\right).
$$
Then, the following chain of equalities holds,
\begin{flushleft}
$
\mathrm{sc}^{([1],2)}_{\boldsymbol{\mathcal{B}}^{(2)}, \varphi(s)}\left(
f^{(2)\flat}_{s}\left(
\mathfrak{P}^{(2)}
\right)
\right)
$
\begin{align*}
&=
\mathrm{sc}^{([1],2)}_{\boldsymbol{\mathcal{B}}^{(2)}, \varphi(s)}\left(
f^{(2)\flat}_{s}\left(
\mathfrak{P}^{(2)}_{1}
\right)
\circ^{0\mathbf{Pth}_{\boldsymbol{\mathcal{B}}^{(2)}}}_{\varphi(s)}
f^{(2)\flat}_{s}\left(
\mathfrak{P}^{(2)}_{0}
\right)
\right)
\tag{1}
\\
&=
\mathrm{sc}^{([1],2)}_{\boldsymbol{\mathcal{B}}^{(2)}, \varphi(s)}\left(
f^{(2)\flat}_{s}\left(
\mathfrak{P}^{(2)}_{1}
\right)
\right)
\circ^{0[\mathbf{PT}_{\boldsymbol{\mathcal{B}}^{(1)}}]}_{\varphi(s)}
\mathrm{sc}^{([1],2)}_{\boldsymbol{\mathcal{B}}^{(2)}, \varphi(s)}\left(
f^{(2)\flat}_{s}\left(
\mathfrak{P}^{(2)}_{0}
\right)
\right)
\tag{2}
\\
&=
f^{[1]@}_{s}\left(
\mathrm{sc}^{([1],2)}_{\boldsymbol{\mathcal{A}}^{(2)}, s}\left(
\mathfrak{P}^{(2)}_{1}
\right)
\right)
\circ^{0[\mathbf{PT}_{\boldsymbol{\mathcal{B}}^{(1)}}]}_{\varphi(s)}
f^{[1]@}_{s}\left(
\mathrm{sc}^{([1],2)}_{\boldsymbol{\mathcal{A}}^{(2)}, s}\left(
\mathfrak{P}^{(2)}_{0}
\right)
\right)
\tag{3}
\\
&=
f^{[1]@}_{s}\left(
\mathrm{sc}^{([1],2)}_{\boldsymbol{\mathcal{A}}^{(2)}, s}\left(
\mathfrak{P}^{(2)}_{1}
\right)
\right)
\circ^{0[\mathbf{PT}_{\boldsymbol{\mathcal{B}}^{(1)}}^{\mathbf{f}^{(1)}}]}_{s}
f^{[1]@}_{s}\left(
\mathrm{sc}^{([1],2)}_{\boldsymbol{\mathcal{A}}^{(2)}, s}\left(
\mathfrak{P}^{(2)}_{0}
\right)
\right)
\tag{4}
\\
&=
f^{[1]@}_{s}\left(
\mathrm{sc}^{([1],2)}_{\boldsymbol{\mathcal{A}}^{(2)}, s}\left(
\mathfrak{P}^{(2)}_{1}
\right)
\circ^{0[\mathbf{PT}_{\boldsymbol{\mathcal{A}}^{(1)}}]}_{s}
\mathrm{sc}^{([1],2)}_{\boldsymbol{\mathcal{A}}^{(2)}, s}\left(
\mathfrak{P}^{(2)}_{0}
\right)
\right)
\tag{5}
\\
&=
f^{[1]@}_{s}\left(
\mathrm{sc}^{([1],2)}_{\boldsymbol{\mathcal{A}}^{(2)}, s}\left(
\mathfrak{P}^{(2)}_{1}
\circ^{0\mathbf{Pth}_{\boldsymbol{\mathcal{A}}^{(2)}}}_{s}
\mathfrak{P}^{(2)}_{0}
\right)
\right)
\tag{6}
\\
&=
f^{[1]@}_{s}\left(
\mathrm{sc}^{([1],2)}_{\boldsymbol{\mathcal{A}}^{(2)}, s}\left(
\mathfrak{P}^{(2)}
\right)
\right)
\tag{7}
\end{align*}
\end{flushleft}

The first equality unravels the definition of $f^{(2)\flat}_{s}(\mathfrak{P}^{(2)})$;
the second equality follows from Claim~\ref{CDPthCatAlgCompZ};
the third equality follows by Artinian induction;
the fourth equality recovers the interpretation of teh operation symbol $\circ^{0}_{s}$ in the partial $\Sigma^{\boldsymbol{\mathcal{A}}^{(1)}}$-algebra $[\mathbf{PT}_{\boldsymbol{\mathcal{B}}^{(1)}}^{\mathbf{f}^{(1)}}]$, introduced in Proposition~\ref{PDQPTBDCatAlg};
the fifth equality follows from Claim~\ref{CDPthCatAlgCompZ};
the sixth equality follows from the fact that, according to Proposition~\ref{PDUCatHom}, $\mathrm{sc}^{([1],2)}_{\boldsymbol{\mathcal{A}}^{(1)}}$ is a $\Sigma^{\boldsymbol{\mathcal{A}}^{(1)}}$-homomorphism from $\mathbf{Pth}_{\boldsymbol{\mathcal{A}}^{(2)}}^{(2,1)}$ to $[\mathbf{PT}_{\boldsymbol{\mathcal{A}}^{(1)}}]$;
finally, the last equality recovers the definition of the second-order path $\mathfrak{P}^{(2)}$.

This completes the proof of Item~(1).

\textsf{(2)}
$\mathrm{tg}^{([1],2)}_{\boldsymbol{\mathcal{B}}^{(2)},\varphi}\circ f^{(2)\flat}=f^{[1]@}\circ \mathrm{tg}^{([1],2)}_{\boldsymbol{\mathcal{A}}^{(1)}}$.

This proof is similar to that of item (1).

\textsf{(3)}
$f^{(2)\flat} \circ \mathrm{ip}_{\boldsymbol{\mathcal{A}}^{(2)}}^{(2, [1])\sharp} = \mathrm{ip}_{\boldsymbol{\mathcal{B}}^{(2)}, \varphi}^{(2, [1])\sharp} \circ f^{[1]@}$.

This property follows directly from the definition of the mapping $f^{(2)\flat}$.

\textsf{(4)}
$f^{(2)\flat} \circ \mathrm{ech}^{(2,\mathcal{A}^{(2)})}_{\boldsymbol{\mathcal{A}}^{(2)}} = f^{(2)}$.

This property follows directly from the definition of the mapping $f^{(2)\flat}$.

This completes the proof.

\end{proof}

We next prove that, for a second-order morphism, its second-order path extension mapping introduced in Proposition~\ref{PDPthExt} is a $\Sigma$-homomorphism.

\begin{proposition}
\label{PDPthExtHom}
Let $\mathbf{f}^{(2)}=(\varphi, c, (f^{(i)})_{i\in 3})$ be a second-order morphism from $\boldsymbol{\mathcal{A}}^{(2)}$ to $\boldsymbol{\mathcal{B}}^{(2)}$. Then the second-order path extension mapping of $f^{(2)}$, $f^{(2)\flat} \colon \mathrm{Pth}_{\boldsymbol{\mathcal{A}}^{(2)}} \mor \mathrm{Pth}_{\boldsymbol{\mathcal{B}}^{(2)}, \varphi}$ introduced in Proposition~\ref{PDPthExt} is a $\Sigma$-homomorphism
$$
f^{(2)\flat} \colon \mathbf{Pth}_{\boldsymbol{\mathcal{A}}^{(2)}}^{(0,2)} \mor \mathbf{Pth}_{\boldsymbol{\mathcal{B}}^{(2)}, \varphi}^{\mathbf{f}^{(2)}(0,2)}.
$$
\end{proposition}

\begin{proof}
Let $(\mathbf{s},s)$ be an element of $S^{\star}\times S$ and $\sigma$ an operation symbol in $\Sigma_{\mathbf{s},s}$.

If $\mathbf{s}=\lambda$, then the following chain of equalities holds
\allowdisplaybreaks
\begin{align*}
f_{s}^{(2)\flat}\left(
\sigma^{\mathbf{Pth}_{\boldsymbol{\mathcal{A}}^{(2)}}^{(0,2)}}
\right)
&=
f_{s}^{(2)\flat}\left(
\mathrm{ip}^{(2,[1])\sharp}_{\boldsymbol{\mathcal{A}}^{(2)},s}\left(
\left[
\sigma^{\mathbf{PT}_{\boldsymbol{\mathcal{A}}^{(1)}}}
\right]_{s}
\right)
\right)
\tag{1}
\\
&=
\mathrm{ip}^{(2,[1])\sharp}_{\boldsymbol{\mathcal{B}}^{(2)},\varphi(s)}\left(
f_{s}^{[1]@}\left(
\left[
\sigma^{\mathbf{PT}_{\boldsymbol{\mathcal{A}}^{(1)}}}
\right]_{s}
\right)
\right)
\tag{2}
\\
&=
\mathrm{ip}_{\boldsymbol{\mathcal{B}}^{(1)},\varphi(s)}^{(2,[1])\sharp}\left(
\left[
\sigma^{\mathbf{PT}^{\mathbf{f}^{(1)}}_{\boldsymbol{\mathcal{B}}^{(1)}}}
\right]_{\varphi(s)}
\right)
\tag{3}
\\
&=
\sigma^{\mathbf{Pth}_{\boldsymbol{\mathcal{B}}^{(2)}}^{\mathbf{f}^{(2)}(0,2)}}.
\tag{4}
\end{align*}

The first equality follows from the fact that, for a constant $\sigma$ in $\Sigma_{\lambda, s}$, then, following Remark~\ref{RDConsSigma}, the interpretation of $\sigma$ in $\mathbf{Pth}_{\boldsymbol{\mathcal{A}}^{(2)}}$ is given by $\Sigma^{\mathbf{Pth}_{\boldsymbol{\mathcal{A}}^{(2)}}} = \mathrm{ip}_{s}^{(2,[1])\sharp}([\sigma^{\mathbf{PT}_{\boldsymbol{\mathcal{A}}^{(1)}}}]_{s})$;
the second equality unravels the definition of $f^{(2)\flat}$ at a $(2,[1])$-identity path;
The third equality follows from the fact that, according to Definition~\ref{DQPthExt}, $f^{[1]@}$ is a $\Sigma$-homomorphism from $\mathbf{PT}_{\boldsymbol{\mathcal{A}}^{(1)}}$ to $\mathbf{PT}_{\boldsymbol{\mathcal{B}}^{(1)}}^{\mathbf{f}^{(1)}}$;
finally, the last equality from the fact that, according to Proposition~\ref{PIpDUBCatHom}, $\mathrm{ip}_{\boldsymbol{\mathcal{B}}^{(1)},\varphi(s)}^{(2,[1])\sharp}$ is a $\Sigma^{\boldsymbol{\mathcal{A}}^{(1)}}$-homomorphism from $[\mathbf{PT}^{\mathbf{f}^{(1)}}_{\boldsymbol{\mathcal{B}}^{(1)}}]$ to $\mathbf{Pth}_{\boldsymbol{\mathcal{B}}^{(2)}}^{\mathbf{f}^{(2)}(1,2)}$.

We now consider the case in which $\mathbf{s}\neq\lambda$. Let $(\mathfrak{P}_{j}^{(2)})_{j\in\bb{\mathbf{s}}}$ be a family of second-order paths in $\mathrm{Pth}_{\boldsymbol{\mathcal{A}}^{(2)},\mathbf{s}}$. We consider different cases according to the nature of the family $(\mathfrak{P}_{j}^{(2)})_{j\in\bb{\mathbf{s}}}$. It could be the case that either~(1), for every $j\in\bb{\mathbf{s}}$, $\mathfrak{P}_{j}^{(2)}$ is a $(2,[1])$-identity second-order path or~(2), there exists an index $j\in\bb{\mathbf{s}}$ for which $\mathfrak{P}_{j}^{(2)}$ is not an $(2,[1])$-identity second-order path.

If~(1), then for every $j\in\bb{\mathbf{s}}$, $\mathfrak{P}_{j}^{(2)}$ is equal to $\mathrm{ip}_{\boldsymbol{\mathcal{A}}^{(2)},s_{j}}^{(2,[1])\sharp}([P_{j}]_{s_{j}})$ for some path term class $[P_{j}]_{s_{j}}$ in $[\mathrm{PT}_{\boldsymbol{\mathcal{A}}^{(1)}}]_{s_{j}}$. In this case, the following chain of equalities holds
\begin{flushleft}
$
f_{s}^{(2)\flat}\left(
\sigma^{\mathbf{Pth}_{\boldsymbol{\mathcal{A}}^{(2)}}^{(0,2)}}\left(\left(
\mathfrak{P}_{j}^{(2)}
\right)_{j\in\bb{\mathbf{s}}}\right)
\right)
$
\allowdisplaybreaks
\begin{align*}
&=
f_{s}^{(2)\flat}\left(
\sigma^{\mathbf{Pth}_{\boldsymbol{\mathcal{A}}^{(2)}}^{(0,2)}}\left(\left(
\mathrm{ip}_{\boldsymbol{\mathcal{A}}^{(2)},s_{j}}^{(2,[1])\sharp}\left(
\left[
P_{j}
\right]_{s_{j}}
\right)
\right)_{j\in\bb{\mathbf{s}}}\right)
\right)
\tag{1}
\\
&=
f_{s}^{(2)\flat}\left(
\mathrm{ip}_{\boldsymbol{\mathcal{A}}^{(1)},s}^{(2,[1])\sharp}\left(
\sigma^{\left[\mathbf{PT}_{\boldsymbol{\mathcal{A}}^{(1)}}\right]}\left(\left(
\left[
P_{j}
\right]_{s_{j}}
\right)_{j\in\bb{\mathbf{s}}}\right)
\right)
\right)
\tag{2}
\\
&=
\mathrm{ip}_{\boldsymbol{\mathcal{B}}^{(1)},\varphi(s)}^{(2,[1])\sharp}\left(
f_{s}^{[1]@}\left(
\sigma^{[\mathbf{PT}_{\boldsymbol{\mathcal{A}}^{(1)}}]}\left(\left(
\left[
P_{j}
\right]_{s_{j}}
\right)_{j\in\bb{\mathbf{s}}}\right)
\right)
\right)
\tag{3}
\\
&=
\mathrm{ip}_{\boldsymbol{\mathcal{B}}^{(2)},\varphi(s)}^{(2,[1])\sharp}\left(
\sigma^{\left[\mathbf{PT}_{\boldsymbol{\mathcal{B}}^{(1)}}^{\mathbf{f}^{(1)}}\right]}\left(\left(
f_{s_{j}}^{[1]@}\left(
\left[
P_{j}
\right]_{s_{j}}
\right)\right)_{j\in\bb{\mathbf{s}}}
\right)
\right)
\tag{4}
\\
&=
\sigma^{\mathbf{Pth}_{\boldsymbol{\mathcal{B}}^{(2)}}^{\mathbf{f}^{(2)}(0,2)}}\left(\left(
\mathrm{ip}_{\boldsymbol{\mathcal{B}}^{(2)},\varphi(s_{j})}^{(2,[1])\sharp}\left(
f_{s_{j}}^{[1]@}\left(
\left[
P_{j}
\right]_{s_{j}}
\right)\right)
\right)_{j\in\bb{\mathbf{s}}}\right)
\tag{5}
\\
&=
\sigma^{\mathbf{Pth}_{\boldsymbol{\mathcal{B}}^{(2)}}^{\mathbf{f}^{(2)}(0,2)}}\left(\left(
{f}_{s_{j}}^{(2)\flat}\left(
\mathrm{ip}_{\boldsymbol{\mathcal{A}}^{(2)},s_{j}}^{(2,[1])\sharp}\left(
\left[
P_{j}
\right]_{s_{j}}
\right)
\right)
\right)_{j\in\bb{\mathbf{s}}}\right)
\tag{6}
\\
&=
\sigma^{\mathbf{Pth}_{\boldsymbol{\mathcal{B}}^{(2)}}^{\mathbf{f}^{(2)}(0,2)}}\left(\left(
{f}_{s_{j}}^{(2)\flat}\left(
\mathfrak{P}_{j}^{(2)}
\right)
\right)_{j\in\bb{\mathbf{s}}}\right).
\tag{7}
\end{align*}
\end{flushleft}

The first equality follows from the fact that, for every $j\in\bb{\mathbf{s}}$, $\mathfrak{P}_{j}^{(2)}$ is a $(2,[1])$-identity path;
the second equality follows from the fact that, by Proposition~\ref{PDUSigma}, $\mathrm{ip}_{\boldsymbol{\mathcal{A}}^{(1)}}^{(2,[1])\sharp}$ is a $\Sigma$-homomorphism from $[\mathbf{PT}_{\boldsymbol{\mathcal{A}}^{(1)}}]$ to $\mathbf{Pth}_{\boldsymbol{\mathcal{A}}^{(2)}}^{(0,2)}$;
the third equality unravels the definition of $f^{(2)\flat}$ at a $(2,[1])$-identity path;
the fourth equality follows from the fact that, by Definition~\ref{DQPthExt}, $f^{[1]@}$ is a $\Sigma^{\boldsymbol{\mathcal{A}}^{(1)}}$-homomorphism from $[\mathbf{PT}_{\boldsymbol{\mathcal{A}}^{(1)}}]$ to $[\mathbf{PT}_{\boldsymbol{\mathcal{B}}^{(1)}}^{\mathbf{f}^{(1)}}]$;
the fifth equality follows from the fact that, by Proposition~\ref{PIpDUBCatHom}, $\mathrm{ip}^{(2,[1])\sharp}_{\boldsymbol{\mathcal{B}}^{(2)}, \varphi}$ is a $\Sigma$-homomorphism from $[\mathbf{PT}^{\mathbf{f}^{(1)}}_{\boldsymbol{\mathcal{B}}^{(1)}}]
$ to $\mathbf{Pth}_{\boldsymbol{\mathcal{B}}^{(2)}}^{\mathbf{f}^{(2)}(1,2)}$;
the sixth equality recovers the definition of $f^{(2)\flat}$ at a $(2,[1])$-identity path;
finally, the last equality recovers, for every $j\in\bb{\mathbf{s}}$, the definition of $\mathfrak{P}_{j}$ as a $(2,[1])$-identity path.
This proves Case~(1).

If~(2), i.e., if there exists some index $j\in\bb{\mathbf{s}}$ for which $\mathfrak{P}_{j}^{(2)}$ is not an $(2,[1])$-identity path then, according to Corollary~\ref{CDPthWB}, $\sigma^{\mathbf{Pth}_{\boldsymbol{\mathcal{A}}^{(2)}}^{(0,2)}}((\mathfrak{P}_{j}^{(2)})_{j\in\bb{\mathbf{s}}})$ is an echelonless, head-constant and coherent second-order path. Moreover, according to Proposition~\ref{PDRecov}, the second-order path extraction algorithm from Lemma~\ref{LDPthExtract}, applied to it retrieves the original family $(\mathfrak{P}_{j}^{(2)})_{j\in\bb{\mathbf{s}}}$. Then, by definition of the $S$-sorted mapping $f^{(2)\flat}$,
$$
f^{(2)\flat}_{s}\left(
\sigma^{\mathbf{Pth}_{\boldsymbol{\mathcal{A}}^{(2)}}^{(1,2)}}\left(\left(
\mathfrak{P}_{j}^{(2)}
\right)_{j\in\bb{\mathbf{s}}}\right)
\right)
=
\sigma^{\mathbf{Pth}_{\boldsymbol{\mathcal{B}}^{(2)}}^{\mathbf{f}^{(2)}(0,2)}}\left(\left(
f^{(2)\flat}_{s_{j}}\left(
\mathfrak{P}_{j}^{(2)}
\right)
\right)_{j\in\bb{\mathbf{s}}}\right)
$$

Case~(2) follows.
\end{proof}

We now show the relation between the second-order path extension mapping of a second-order morphism and the mappings $\mathrm{sc}^{(0,2)}$, $\mathrm{tg}^{(0,2)}$ and $\mathrm{ip}^{(2,0)}\sharp$.

\begin{proposition}
\label{PDPthExtDZScTg}
Let $\mathbf{f}^{(2)}=(\varphi, c, (f^{(i)})_{i\in 3})$ be a second-order morphism from $\boldsymbol{\mathcal{A}}^{(2)}$ to $\boldsymbol{\mathcal{B}}^{(2)}$. Then the following equalities holds
\begin{align*}
\mathrm{sc}_{\boldsymbol{\mathcal{B}}^{(2)},\varphi}^{(0,2)} \circ f^{(2)\flat}
&=
f^{(0)\sharp}\circ\mathrm{sc}_{\boldsymbol{\mathcal{A}}^{(2)}}^{(0,2)}
&&\mbox{and}&
\mathrm{tg}_{\boldsymbol{\mathcal{B}}^{(2)},\varphi}^{(0,2)} \circ f^{(2)\flat}
&=
f^{(0)\sharp}\circ\mathrm{tg}_{\boldsymbol{\mathcal{A}}^{(2)}}^{(0,2)}.
\end{align*}
\end{proposition}

\begin{proof}
The following chain of equalities holds
\allowdisplaybreaks
\begin{align*}
\mathrm{sc}_{\boldsymbol{\mathcal{B}}^{(2)},\varphi}^{(0,2)}
\circ
f^{(2)\flat}
&=
\mathrm{sc}^{(0,[1])}_{\boldsymbol{\mathcal{B}}^{(2)},\varphi}
\circ
\mathrm{ip}^{([1],X)@}_{\boldsymbol{\mathcal{B}}^{(2)},\varphi}
\circ
\mathrm{sc}^{([1],2)}_{\boldsymbol{\mathcal{B}}^{(2)},\varphi}
\circ
f^{(2)\flat}
\tag{1}
\\
&=
\mathrm{sc}^{(0,[1])}_{\boldsymbol{\mathcal{B}}^{(2)},\varphi}
\circ
\mathrm{ip}^{([1],X)@}_{\boldsymbol{\mathcal{B}}^{(2)},\varphi}
\circ
f^{[1]@}
\circ
\mathrm{sc}^{([1],2)}_{\boldsymbol{\mathcal{A}}^{(2)}}
\tag{2}
\\
&=
\mathrm{sc}^{(0,[1])}_{\boldsymbol{\mathcal{B}}^{(2)},\varphi}
\circ
f^{[1]@}
\circ
\mathrm{ip}^{([1],X)@}_{\boldsymbol{\mathcal{A}}^{(2)}}
\circ
\mathrm{sc}^{([1],2)}_{\boldsymbol{\mathcal{A}}^{(2)}}
\tag{3}
\\
&=
f^{(0)\sharp}
\circ
\mathrm{sc}^{(0,[1])}_{\boldsymbol{\mathcal{A}}^{(2)}}
\circ
\mathrm{ip}^{([1],X)@}_{\boldsymbol{\mathcal{A}}^{(2)}}
\circ
\mathrm{sc}^{([1],2)}_{\boldsymbol{\mathcal{A}}^{(2)}}
\tag{4}
\\
&=
f^{(0)\sharp}
\circ
\mathrm{sc}_{\boldsymbol{\mathcal{A}}^{(2)}}^{(0,2)}
\tag{5}
\end{align*}

The first equality unravels the definition of the mapping \(\mathrm{sc}^{(0,2)}_{\boldsymbol{\mathcal{B}}^{(2)}, \varphi}\);
the second equality follows from Proposition~\ref{PDPthExt};
the third equality follows from Proposition~\ref{PQPthExtQIp};
the fourth equality follows from Proposition~\ref{PQPthExtScTg};
finally, the last equality recovers the definition of the mapping \(\mathrm{sc}_{\boldsymbol{\mathcal{A}}^{(2)}}^{(0,2)}\).

A similar argument applies to the $(0,2)$-target mapping.
\end{proof}

\begin{proposition}
\label{PDPthExtDZIp}
Let $\mathbf{f}^{(2)}=(\varphi, c, (f^{(i)})_{i\in 3})$ be a second-order morphism from $\boldsymbol{\mathcal{A}}^{(2)}$ to $\boldsymbol{\mathcal{B}}^{(2)}$. Then $f^{(2)\flat} \circ \mathrm{ip}_{\boldsymbol{\mathcal{A}}^{(2)}}^{(2,0)\sharp} = \mathrm{ip}_{\boldsymbol{\mathcal{B}}^{(2)}, \varphi}^{(2,0)\sharp} \circ f^{(0)\sharp}$.
\end{proposition}

\begin{proof}
The following chain of equalities holds
\allowdisplaybreaks
\begin{align*}
f^{(2)\flat}
\circ
\mathrm{ip}_{\boldsymbol{\mathcal{B}}^{(2)},\varphi}^{(2,0)\sharp}
&=
f^{(2)\flat}
\circ
\mathrm{ip}^{(2,[1])\sharp}_{\boldsymbol{\mathcal{B}}^{(2)},\varphi}
\circ
\mathrm{CH}^{[1]}_{\boldsymbol{\mathcal{B}}^{(2)},\varphi}
\circ
\mathrm{ip}^{([1], 0)\sharp}_{\boldsymbol{\mathcal{B}}^{(2)},\varphi}
\tag{1}
\\
&=
\mathrm{ip}^{(2,[1])\sharp}_{\boldsymbol{\mathcal{B}}^{(2)},\varphi}
\circ
f^{[1]@}
\circ
\mathrm{CH}^{[1]}_{\boldsymbol{\mathcal{B}}^{(2)},\varphi}
\circ
\mathrm{ip}^{([1], 0)\sharp}_{\boldsymbol{\mathcal{B}}^{(2)},\varphi}
\tag{2}
\\
&=
\mathrm{ip}^{(2,[1])\sharp}_{\boldsymbol{\mathcal{B}}^{(2)},\varphi}
\circ
\mathrm{CH}^{[1]}_{\boldsymbol{\mathcal{B}}^{(2)},\varphi}
\circ
f^{[1]@}
\circ
\mathrm{ip}^{([1], 0)\sharp}_{\boldsymbol{\mathcal{B}}^{(2)},\varphi}
\tag{3}
\\
&=
\mathrm{ip}^{(2,[1])\sharp}_{\boldsymbol{\mathcal{B}}^{(2)},\varphi}
\circ
\mathrm{CH}^{[1]}_{\boldsymbol{\mathcal{B}}^{(2)},\varphi}
\circ
\mathrm{ip}^{([1], 0)\sharp}_{\boldsymbol{\mathcal{B}}^{(2)},\varphi}
\circ
f^{(0)\sharp}
\tag{4}
\\
&=
\mathrm{ip}_{\boldsymbol{\mathcal{A}}^{(2)}}^{(2,0)\sharp}
\circ
f^{(0)\sharp}
\tag{5}
\end{align*}

The first equality unravels the definition of the mapping \(\mathrm{ip}^{(2,0)\sharp}_{\boldsymbol{\mathcal{B}}^{(2)}, \varphi}\);
the second equality follows from Proposition~\ref{PDPthExt};
the third equality follows from Proposition~\ref{PQPthExtQCH};
the fourth equality follows from Proposition~\ref{PQPthExtIp};
finally, the last equality recovers the definition of the mapping \(\mathrm{ip}_{\boldsymbol{\mathcal{A}}^{(2)}}^{(2,0)\sharp}\).
\end{proof}

Finally, given a second-order morphism $\mathbf{f}^{(2)}$ from $\boldsymbol{\mathcal{A}}^{(2)}$ to $\boldsymbol{\mathcal{B}}^{(2)}$, we equip the $S$-sorted sets $\mathrm{Pth}_{\boldsymbol{\mathcal{B}}^{(2)}, \varphi}$ and $\llbracket\mathrm{Pth}_{\boldsymbol{\mathcal{B}}^{(2)}}\rrbracket_{\varphi}$ with a structure of partial $\Sigma^{\boldsymbol{\mathcal{A}}^{(1)}}$-algebra and show that the second-order path extension mapping is a $\Sigma^{\boldsymbol{\mathcal{A}}^{(1)}}$-homomorphism.

\begin{proposition}\label{PDPthBCatAlg}
Let $\mathbf{f}^{(2)}=(\varphi, c, (f^{(i)})_{i\in 3})$ be a second-order morphism from $\boldsymbol{\mathcal{A}}^{(2)}$ to $\boldsymbol{\mathcal{B}}^{(2)}$. Then the $S$-sorted set $\mathrm{Pth}_{\boldsymbol{\mathcal{B}}^{(2)}, \varphi}$ is equipped, in a natural way, with a structure of partial $\Sigma^{\boldsymbol{\mathcal{A}}^{(1)}}$-algebra.
\end{proposition}

\begin{proof}
Let us denote by $\mathbf{Pth}_{\boldsymbol{\mathcal{B}}^{(2)}}^{\mathbf{f}^{(2)}(1,2)}$ the partial $\Sigma^{\boldsymbol{\mathcal{A}}^{(1)}}$-algebra defined as follows 

\textsf{(1)}
The underlying $S$-sorted set of $\mathbf{Pth}_{\boldsymbol{\mathcal{B}}^{(2)}}^{\mathbf{f}^{(2)}(1,2)}$ is $\mathrm{Pth}_{\boldsymbol{\mathcal{B}}^{(2)}, \varphi}$.

\textsf{(2)}
For every $(\mathbf{s},s) \in S^{\star}\times S$ and every operation symbol $\sigma \in \Sigma_{\mathbf{s},s}$, the operation $\sigma^{\mathbf{Pth}_{\boldsymbol{\mathcal{B}}^{(2)}}^{\mathbf{f}^{(2)}(1,2)}}$ is given by the interpretation of $\sigma$ in the $\Sigma$-algebra $\mathbf{Pth}_{\boldsymbol{\mathcal{B}}^{(2)}}^{\mathbf{f}^{(2)}(0,2)}$ introduced in Remark~\ref{RDDSigmaAlg}.

\textsf{(3)}
For every $s \in S$ and every $\mathfrak{p} \in \mathcal{A}^{(1)}_{s}$, the constant $\mathfrak{p}^{\mathbf{Pth}_{\boldsymbol{\mathcal{B}}^{(2)}}^{\mathbf{f}^{(2)}(1,2)}}$ si given by
$$
\mathfrak{p}^{\mathbf{Pth}_{\boldsymbol{\mathcal{B}}^{(2)}}^{\mathbf{f}^{(2)}(1,2)}}
=
f^{(2)\flat}_{s}\left(
\mathfrak{p}^{\mathbf{Pth}_{\boldsymbol{\mathcal{A}}^{(2)}}}
\right).
$$

\textsf{(4)}
For every $s \in S$, the interpretations fo the operations $\mathrm{sc}_{s}^{0}$ and $\mathrm{tg}_{s}^{0}$ are given by
\begin{align*}
\mathrm{sc}_{s}^{0\mathbf{Pth}_{\boldsymbol{\mathcal{B}}^{(2)}}^{\mathbf{f}^{(2)}(1,2)}}
&=
\mathrm{sc}_{\varphi(s)}^{0\mathbf{Pth}_{\boldsymbol{\mathcal{B}}^{(2)}}}
&&\mbox{and}&
\mathrm{tg}_{s}^{0\mathbf{Pth}_{\boldsymbol{\mathcal{B}}^{(2)}}^{\mathbf{f}^{(2)}(1,2)}}
&=
\mathrm{tg}_{\varphi(s)}^{0\mathbf{Pth}_{\boldsymbol{\mathcal{B}}^{(2)}}}
\end{align*}

\textsf{(5)}
Similarly, the interpretation of the partial binary operation $\circ_{s}^{0}$ is given by 
$$
\circ_{s}^{0\mathbf{Pth}_{\boldsymbol{\mathcal{B}}^{(2)}}^{\mathbf{f}^{(2)}(1,2)}}
=
\circ_{\varphi(s)}^{0\mathbf{Pth}_{\boldsymbol{\mathcal{B}}^{(2)}}}.
$$
\end{proof}

\begin{proposition}\label{PDQPthBCatAlg}
Let $\mathbf{f}^{(2)}=(\varphi, c, (f^{(i)})_{i\in 3})$ be a second-order morphism from $\boldsymbol{\mathcal{A}}^{(2)}$ to $\boldsymbol{\mathcal{B}}^{(2)}$. Then the $S$-sorted set $\llbracket\mathrm{Pth}_{\boldsymbol{\mathcal{B}}^{(2)}}\rrbracket_{\varphi}$ is equipped, in a natural way, with a structure of partial $\Sigma^{\boldsymbol{\mathcal{A}}^{(1)}}$-algebra.
\end{proposition}

\begin{proof}
Let us denote by $\llbracket\mathbf{Pth}_{\boldsymbol{\mathcal{B}}^{(2)}}^{\mathbf{f}^{(2)}(1,2)}\rrbracket$ the partial $\Sigma^{\boldsymbol{\mathcal{A}}^{(1)}}$-algebra defined as follows 

\textsf{(1)}
The underlying $S$-sorted set of $\llbracket\mathbf{Pth}_{\boldsymbol{\mathcal{B}}^{(2)}}^{\mathbf{f}^{(2)}(1,2)}\rrbracket$ is $\llbracket\mathrm{Pth}_{\boldsymbol{\mathcal{B}}^{(2)}}\rrbracket_{\varphi}$.

\textsf{(2)}
For every $(\mathbf{s},s) \in S^{\star}\times S$ and every operation symbol $\sigma \in \Sigma_{\mathbf{s},s}$, the operation $\sigma^{\llbracket\mathbf{Pth}_{\boldsymbol{\mathcal{B}}^{(2)}}^{\mathbf{f}^{(2)}(1,2)}\rrbracket}$ is given by the interpretation of $\sigma$ in the $\Sigma$-algebra $\llbracket\mathbf{Pth}_{\boldsymbol{\mathcal{B}}^{(2)}}^{\mathbf{f}^{(2)}(0,2)}\rrbracket$ introduced in Remark~\ref{RDDSigmaAlg}.

\textsf{(3)}
For every $s \in S$ and every $\mathfrak{p} \in \mathcal{A}^{(1)}_{s}$, the constant $\mathfrak{p}^{\llbracket\mathbf{Pth}_{\boldsymbol{\mathcal{B}}^{(2)}}^{\mathbf{f}^{(2)}(1,2)}\rrbracket}$ si given by
$$
\mathfrak{p}^{\llbracket\mathbf{Pth}_{\boldsymbol{\mathcal{B}}^{(2)}}^{\mathbf{f}^{(2)}(1,2)}\rrbracket}
=
\left\llbracket
f^{(2)\flat}_{s}\left(
\mathfrak{p}^{\mathbf{Pth}_{\boldsymbol{\mathcal{A}}^{(2)}}}
\right)
\right\rrbracket_{\varphi(s)}
$$

\textsf{(4)}
For every $s \in S$, the interpretations fo the operations $\mathrm{sc}_{s}^{0}$ and $\mathrm{tg}_{s}^{0}$ are given by
\begin{align*}
\mathrm{sc}_{s}^{0\llbracket\mathbf{Pth}_{\boldsymbol{\mathcal{B}}^{(2)}}^{\mathbf{f}^{(2)}(1,2)}\rrbracket}
&=
\mathrm{sc}_{\varphi(s)}^{0\llbracket\mathbf{Pth}_{\boldsymbol{\mathcal{B}}^{(2)}}\rrbracket}
&&\mbox{and}&
\mathrm{tg}_{s}^{0\llbracket\mathbf{Pth}_{\boldsymbol{\mathcal{B}}^{(2)}}^{\mathbf{f}^{(2)}(1,2)}\rrbracket}
&=
\mathrm{tg}_{\varphi(s)}^{0\llbracket\mathbf{Pth}_{\boldsymbol{\mathcal{B}}^{(2)}}\rrbracket}
\end{align*}

\textsf{(5)}
Similarly, the interpretation of the partial binary operation $\circ_{s}^{0}$ is given by 
$$
\circ_{s}^{0\llbracket\mathbf{Pth}_{\boldsymbol{\mathcal{B}}^{(2)}}^{\mathbf{f}^{(2)}(1,2)}\rrbracket}
=
\circ_{\varphi(s)}^{0\llbracket\mathbf{Pth}_{\boldsymbol{\mathcal{B}}^{(2)}}\rrbracket}.
$$
\end{proof}

We now show that the many-sorted projection mapping $\mathrm{pr}^{\llbracket\cdot\rrbracket}_{\boldsymbol{\mathcal{B}}^{(2)}, \varphi}$ and the many-sorted mappings $\mathrm{sc}^{([1],2)}_{\boldsymbol{\mathcal{B}}^{(2)}, \varphi}$, $\mathrm{tg}^{([1],2)}_{\boldsymbol{\mathcal{B}}^{(2)}, \varphi}$ and $\mathrm{ip}^{(2,[1])}_{\boldsymbol{\mathcal{B}}^{(2)}, \varphi}$, of $([1], 2)$-source, $([1], 2)$-target and $(2, [1])$-identity second-order paths, respectively, are $\Sigma^{\boldsymbol{\mathcal{A}}^{(1)}}$ -homomorphisms.

\begin{proposition}\label{PDPrBCatHom}
Let $\mathbf{f}^{(2)}=(\varphi, c, (f^{(i)})_{i\in 3})$ be a second-order morphism from $\boldsymbol{\mathcal{A}}^{(2)}$ to $\boldsymbol{\mathcal{B}}^{(2)}$. Then $\mathrm{pr}^{\llbracket\cdot\rrbracket}_{\boldsymbol{\mathcal{B}}^{(2)}, \varphi}$ is a $\Sigma^{\boldsymbol{\mathcal{A}}^{(1)}}$-homomorphism from  $\mathbf{Pth}_{\boldsymbol{\mathcal{B}}^{(2)}}^{\mathbf{f}^{(2)}(1,2)}$ to $\llbracket\mathbf{Pth}_{\boldsymbol{\mathcal{B}}^{(2)}}^{\mathbf{f}^{(2)}(1,2)}\rrbracket$.
\end{proposition}

\begin{proof}
We prove that $\mathrm{pr}_{\boldsymbol{\mathcal{B}}^{(2)}, \varphi}^{\llbracket \cdot \rrbracket}$ is compatible with every operation symbol in $\Sigma^{\boldsymbol{\mathcal{A}}^{(1)}}$.

{\sffamily The mapping $\mathrm{pr}_{\boldsymbol{\mathcal{B}}^{(2)}, \varphi}^{\llbracket \cdot \rrbracket}$ is a $\Sigma$-homomorphism.}

Note that $\mathrm{pr}_{\boldsymbol{\mathcal{B}}^{(2)}, \varphi}^{\llbracket \cdot \rrbracket} = \mathbf{c}_{\mathfrak{d}}^{\ast} (\mathrm{pr}_{\boldsymbol{\mathcal{B}}^{(2)}}^{\llbracket \cdot \rrbracket})$. By Reamrk~\ref{RDDSigmaAlg}, $\mathrm{pr}_{\boldsymbol{\mathcal{B}}^{(2)}}^{\llbracket \cdot \rrbracket}$ is a $\Lambda$-homomorphism. Therefore, it follows from Proposition~\ref{PFunSig} that the mapping $\mathrm{pr}_{\boldsymbol{\mathcal{B}}^{(2)}, \varphi}^{\llbracket \cdot \rrbracket}$ is a $\Sigma$-homomorphism.

{\sffamily The mapping $\mathrm{pr}_{\boldsymbol{\mathcal{B}}^{(2)}, \varphi}^{\llbracket \cdot \rrbracket}$ is compatible with the first-order rewrite rules.}

Let $s$ be a sort in $S$ and $\mathfrak{p}$ a rewrite rule in $\mathcal{A}_{s}^{(1)}$. Thus,
$$
\left\llbracket
\mathfrak{p}^{\mathbf{Pth}_{\boldsymbol{\mathcal{B}}^{(2)}}^{\mathbf{f}^{(2)}(1,2)}}
\right\rrbracket_{\varphi(s)}
=
\left\llbracket
f_{s}^{(2)\flat}\left(
\mathfrak{p}^{\mathbf{Pth}_{\boldsymbol{\mathcal{B}}^{(2)}}}
\right)
\right\rrbracket_{\varphi(s)}
=
\mathfrak{p}^{\llbracket\mathbf{Pth}_{\boldsymbol{\mathcal{B}}^{(2)}}^{\mathbf{f}^{(2)}(1,2)}\rrbracket}.
$$

Hence, $\mathrm{pr}_{\boldsymbol{\mathcal{B}}^{(2)}, \varphi}^{\llbracket \cdot \rrbracket}$ is compatible with the first-order rewrite rules.

{\sffamily The mapping $\mathrm{pr}_{\boldsymbol{\mathcal{B}}^{(2)}, \varphi}^{\llbracket \cdot \rrbracket}$ is compatible with the $0$-source.}

Let $s$ be a sort in $S$ and let us consider the $0$-source operation symbol $\mathrm{sc}_{s}^{0}$ in $\Sigma^{\boldsymbol{\mathcal{A}}^{(1)}}_{s,s}$. Let $\mathfrak{P}^{(2)}$ be a second-order path in $\mathrm{Pth}_{\boldsymbol{\mathcal{B}}^{(2)}, \varphi(s)}$.

The following chain of equalities holds
\allowdisplaybreaks
\begin{align*}
\left\llbracket
\mathrm{sc}_{s}^{0\mathbf{Pth}_{\boldsymbol{\mathcal{B}}^{(2)}}^{\mathbf{f}^{(2)}(1,2)}}\left(
\mathfrak{P}^{(2)}
\right)
\right\rrbracket_{\varphi(s)}
&=
\left\llbracket
\mathrm{sc}_{\varphi(s)}^{0\mathbf{Pth}_{\boldsymbol{\mathcal{B}}^{(2)}}}\left(
\mathfrak{P}^{(2)}
\right)
\right\rrbracket_{\varphi(s)}
\tag{1}
\\
&=
\mathrm{sc}_{\varphi(s)}^{0\llbracket\mathbf{Pth}_{\boldsymbol{\mathcal{B}}^{(2)}}\rrbracket}\left(
\left\llbracket
\mathfrak{P}^{(2)}
\right\rrbracket_{\varphi(s)}
\right)
\tag{2}
\\
&=
\mathrm{sc}_{s}^{0\llbracket\mathbf{Pth}_{\boldsymbol{\mathcal{B}}^{(2)}}^{\mathbf{f}^{(2)}(1,2)}\rrbracket}\left(
\left\llbracket
\mathfrak{P}^{(2)}
\right\rrbracket_{\varphi(s)}
\right).
\tag{3}
\end{align*}
The first equality unravels the interpretation of the operation symbol $\mathrm{sc}_{s}^{0}$ in the partial $\Sigma^{\boldsymbol{\mathcal{A}}^{(1)}}$-algebra $\mathbf{Pth}_{\boldsymbol{\mathcal{B}}^{(2)}}^{\mathbf{f}^{(2)}(1,2)}$, introduced in Proposition~\ref{PDPthBCatAlg};
the second equality follows from the fact that, according to Proposition~\ref{PDVDCatAlg}, $\mathrm{pr}^{\llbracket\cdot\rrbracket}_{\boldsymbol{\mathcal{B}}^{(2)}}$ is a $\Lambda^{\boldsymbol{\mathcal{B}}^{(1)}}$-homomorphism from $\mathbf{Pth}_{\boldsymbol{\mathcal{B}}^{(2)}}^{(1,2)}$ to $\llbracket\mathbf{Pth}_{\boldsymbol{\mathcal{B}}^{(2)}}^{(1,2)}\rrbracket$;
finally, the last equality recovers the interpretation of the operation symbol $\mathrm{sc}_{s}^{0}$ in the partial $\Sigma^{\boldsymbol{\mathcal{A}}^{(1)}}$-algebra $\llbracket\mathbf{Pth}_{\boldsymbol{\mathcal{B}}^{(2)}}^{\mathbf{f}^{(2)}(1,2)}\rrbracket$, introduced in Proposition~\ref{PDQPthBCatAlg}.

Hence, $\mathrm{pr}_{\boldsymbol{\mathcal{B}}^{(2)}, \varphi}^{\llbracket \cdot \rrbracket}$ is compatible with the $0$-source operation.

{\sffamily The mapping $\mathrm{pr}_{\boldsymbol{\mathcal{B}}^{(2)}, \varphi}^{\llbracket \cdot \rrbracket}$ is compatible with the $0$-target.}

Let $s$ be a sort in $S$ and let us consider the $0$-target operation symbol $\mathrm{tg}_{s}^{0}$ in $\Sigma^{\boldsymbol{\mathcal{A}}^{(1)}}_{s,s}$. Let $\mathfrak{P}^{(2)}$ be a second-order path in $\mathrm{Pth}_{\boldsymbol{\mathcal{B}}^{(2)}, \varphi(s)}$, then the following equality holds
$$
\left\llbracket
\mathrm{tg}_{s}^{0\mathbf{Pth}_{\boldsymbol{\mathcal{B}}^{(2)}}^{\mathbf{f}^{(2)}(1,2)}}\left(
\mathfrak{P}^{(2)}
\right)
\right\rrbracket_{\varphi(s)}
=
\mathrm{tg}_{s}^{0\llbracket\mathbf{Pth}_{\boldsymbol{\mathcal{B}}^{(2)}}^{\mathbf{f}^{(2)}(1,2)}\rrbracket}\left(
\left\llbracket
\mathfrak{P}^{(2)}
\right\rrbracket_{\varphi(s)}
\right).
$$

The proof of this case is identical to that of the $0$-source.

Hence, $\mathrm{pr}_{\boldsymbol{\mathcal{B}}^{(2)}, \varphi}^{\llbracket \cdot \rrbracket}$ is compatible with the $0$-target operation.

{\sffamily The mapping $\mathrm{pr}_{\boldsymbol{\mathcal{B}}^{(2)}, \varphi}^{\llbracket \cdot \rrbracket}$ is compatible with the $0$-composition.}

Let $s$ be a sort in $S$ and let us consider the $0$-composition operation symbol $\circ_{s}^{0}$ in $\Sigma_{ss,s}^{\boldsymbol{\mathcal{A}}^{(1)}}$. Let $\mathfrak{P}^{(2)}$ and $\mathfrak{Q}^{(2)}$ be two second-order paths in $\mathrm{Pth}_{\boldsymbol{\mathcal{B}}^{(2)}, \varphi(s)}$ such that
$$
\mathrm{sc}_{\boldsymbol{\mathcal{B}}^{(2)}, \varphi(s)}^{(0,2)}\left(\mathfrak{Q}^{(2)}\right)
=
\mathrm{tg}_{\boldsymbol{\mathcal{B}}^{(2)}, \varphi(s)}^{(0,2)}\left(\mathfrak{P}^{(2)}\right).
$$

Then the following equality holds
$$
\left\llbracket
\mathfrak{Q}^{(2)}
\circ_{s}^{0\mathbf{Pth}_{\boldsymbol{\mathcal{B}}^{(2)}}^{\mathbf{f}^{(2)}(1,2)}}
\mathfrak{P}^{(2)}
\right\rrbracket_{\varphi(s)}
=
\left\llbracket
\mathfrak{Q}^{(2)}
\right\rrbracket_{\varphi(s)}
\circ_{s}^{0\llbracket\mathbf{Pth}_{\boldsymbol{\mathcal{B}}^{(2)}}^{\mathbf{f}^{(2)}(1,2)}\rrbracket}
\left\llbracket
\mathfrak{P}^{(2)}
\right\rrbracket_{\varphi(s)}.
$$

The proof of this case is identical to that of the $0$-source.

Hence, $\mathrm{pr}_{\boldsymbol{\mathcal{B}}^{(2)}, \varphi}^{\llbracket \cdot \rrbracket}$ is compatible with the $0$-composition operation.

This completes the proof.
\end{proof}

\begin{proposition}\label{PScTgUDBCatHom}
Let $\mathbf{f}^{(2)}=(\varphi, c, (f^{(i)})_{i\in 3})$ be a second-order morphism from $\boldsymbol{\mathcal{A}}^{(2)}$ to $\boldsymbol{\mathcal{B}}^{(2)}$. Then $\mathrm{sc}^{([1],2)}_{\boldsymbol{\mathcal{B}}^{(2)}, \varphi}$ and $\mathrm{tg}^{([1],2)}_{\boldsymbol{\mathcal{B}}^{(2)}, \varphi}$ are $\Sigma^{\boldsymbol{\mathcal{A}}^{(1)}}$-homomorphism from $\mathbf{Pth}_{\boldsymbol{\mathcal{B}}^{(2)}}^{\mathbf{f}^{(2)}(1,2)}$ to $[\mathbf{PT}_{\boldsymbol{\mathcal{B}}^{(1)}}^{\mathbf{f}^{(1)}}]$.
\end{proposition}

\begin{proof}
We prove that $\mathrm{sc}^{([1],2)}_{\boldsymbol{\mathcal{B}}^{(2)}, \varphi}$ is compatible with every operation symbol in $\Sigma^{\boldsymbol{\mathcal{A}}^{(1)}}$.

{\sffamily The mapping $\mathrm{sc}^{([1],2)}_{\boldsymbol{\mathcal{B}}^{(2)}, \varphi}$ is a $\Sigma$-homomorphism.}

Note that $\mathrm{sc}^{([1],2)}_{\boldsymbol{\mathcal{B}}^{(2)}, \varphi} = \mathbf{c}_{\mathfrak{d}}^{\ast} (\mathrm{sc}^{([1],2)}_{\boldsymbol{\mathcal{B}}^{(2)}})$. By Proposition~\ref{PDVDU}, the mapping $\mathrm{sc}^{([1],2)}_{\boldsymbol{\mathcal{B}}^{(2)}}$ is a $\Lambda^{\boldsymbol{\mathcal{B}}^{(1)}}$-homomorphism, thus in particular a $\Lambda$-homomorphism. Therefore, it follows from Proposition~\ref{PFunSig} that the mapping $\mathrm{sc}^{([1],2)}_{\boldsymbol{\mathcal{B}}^{(2)}, \varphi}$ is a $\Sigma$-homomorphism.

{\sffamily The mapping $\mathrm{sc}^{([1],2)}_{\boldsymbol{\mathcal{B}}^{(2)}, \varphi}$ is compatible with the first-order rewrite rules.}

Let $s$ be a sort in $S$ and $\mathfrak{p}$ a rewrite rule in $\mathcal{A}_{s}^{(1)}$. Thus, the following chain of equalities holds
\begin{flushleft}
$
\mathrm{sc}^{([1],2)}_{\boldsymbol{\mathcal{B}}^{(2)}, \varphi(s)}\left(
\mathfrak{p}^{\mathbf{Pth}_{\boldsymbol{\mathcal{B}}^{(2)}}^{\mathbf{f}^{(2)}(1,2)}}
\right)
$
\allowdisplaybreaks
\begin{align*}
&=
\mathrm{sc}^{([1],2)}_{\boldsymbol{\mathcal{B}}^{(2)}, \varphi(s)}\left(
f_{s}^{(2)\flat}\left(
\mathfrak{p}^{\mathbf{Pth}_{\boldsymbol{\mathcal{A}}^{(2)}}}
\right)
\right)
\tag{1}
\\
&=
f_{s}^{[1]@}\left(
\mathrm{sc}_{\boldsymbol{\mathcal{A}}^{(2)},s}^{([1],2)}\left(
\mathfrak{p}^{\mathbf{Pth}_{\boldsymbol{\mathcal{A}}^{(2)}}}
\right)
\right)
\tag{2}
\\
&=
f_{s}^{[1]@}\left(
\mathfrak{p}^{[\mathbf{PT}_{\boldsymbol{\mathcal{A}}^{(1)}}]}
\right)
\tag{3}
\\
&=
f_{s}^{[1]@}\left(
\left[
\mathfrak{p}^{\mathbf{PT}_{\boldsymbol{\mathcal{A}}^{(1)}}}
\right]_{s}
\right)
\tag{4}
\\
&=
\left[
\mathrm{CH}^{(1)}_{\boldsymbol{\mathcal{B}}^{(1)}, \varphi(s)}\left(
f_{s}^{(1)\flat}\left(
\mathrm{ip}^{(1,X)@}_{\boldsymbol{\mathcal{A}}^{(1)}, s}\left(
\mathfrak{p}^{\mathbf{PT}_{\boldsymbol{\mathcal{A}}^{(1)}}}
\right)
\right)
\right)
\right]_{\varphi(s)}
\tag{5}
\\
&=
\left[
\mathrm{CH}^{(1)}_{\boldsymbol{\mathcal{B}}^{(1)}, \varphi(s)}\left(
f_{s}^{(1)\flat}\left(
\mathfrak{p}^{\mathbf{Pth}_{\boldsymbol{\mathcal{A}}^{(1)}}}
\right)
\right)
\right]_{\varphi(s)}
\tag{6}
\\
&=
\mathfrak{p}^{[\mathbf{PT}_{\boldsymbol{\mathcal{B}}^{(1)}}^{\mathbf{f}^{(1)}}]}.
\tag{7}
\end{align*}
\end{flushleft}

In the just stated chain of equalities, the first equality unravels the definition of the interpretation of the constant symbol $\mathfrak{p}$ in the partial $\Sigma^{\boldsymbol{\mathcal{A}}^{(1)}}$-algebra $\mathbf{Pth}^{\mathbf{f}^{(2)}(1,2)}_{\boldsymbol{\mathcal{B}}^{(2)}}$, introduced in Proposition~\ref{PDPthBCatAlg};
the second equality follows from Proposition~\ref{PDPthExt};
the third equality follows from the fact that, according to Proposition~\ref{PDUCatHom}, $\mathrm{sc}^{([1],2)}_{\boldsymbol{\mathcal{A}}^{(2)}}$ is a $\Sigma^{\boldsymbol{\mathcal{A}}^{(1)}}$-homomorphism from $\mathbf{Pth}_{\boldsymbol{\mathcal{A}}^{(2)}}$ to $[\mathbf{PT}_{\boldsymbol{\mathcal{A}}^{(1)}}]$;
the fourth equality unravels the definition of the interpretation of the constant symbol $\mathfrak{p}$ in the partial $\Sigma^{\boldsymbol{\mathcal{A}}^{(1)}}$-algebra $[\mathbf{PT}_{\boldsymbol{\mathcal{A}}^{(1)}}]$, introduced in Proposition~\ref{PPTQCatAlg};
the fifth equality unravels the definition of the $s$-th component of the mapping $f^{[1]@}$ at a path term class;
the sixth equality follows from the fact that, by Definition~\ref{DIp}, $\mathrm{ip}^{(1,X)@}_{\boldsymbol{\mathcal{A}}^{(1)}}$ is a $\Sigma^{\boldsymbol{\mathcal{A}}^{(1)}}$-homomorphism from $\mathbf{PT}_{\boldsymbol{\mathcal{A}}^{(1)}}$ to $\mathbf{Pth}_{\boldsymbol{\mathcal{A}}^{(1)}}$; 
finally, the last equality  recovers the definition of the interpretation of the constant symbol $\mathfrak{p}$ in the partial $\Sigma^{\boldsymbol{\mathcal{A}}^{(1)}}$-algebra $[\mathbf{PT}^{\mathbf{f}^{(1)}}_{\boldsymbol{\mathcal{B}}^{(1)}}]$, introduced in Proposition~\ref{PPTBCatAlg}.

Hence, $\mathrm{sc}^{([1],2)}_{\boldsymbol{\mathcal{B}}^{(2)}, \varphi}$ is compatible with the first-order rewrite rules.

{\sffamily The mapping $\mathrm{sc}^{([1],2)}_{\boldsymbol{\mathcal{B}}^{(2)}, \varphi}$ is compatible with the $0$-source.}

Let $s$ be a sort in $S$ and let us consider the $0$-source operation symbol $\mathrm{sc}_{s}^{0}$ in $\Sigma^{\boldsymbol{\mathcal{A}}^{(1)}}_{s,s}$. Let $\mathfrak{P}^{(2)}$ be a second-order path in $\mathrm{Pth}_{\boldsymbol{\mathcal{B}}^{(2)}, \varphi(s)}$.

The following chain of equalities holds
\allowdisplaybreaks
\begin{align*}
\mathrm{sc}^{([1],2)}_{\boldsymbol{\mathcal{B}}^{(2)}, \varphi(s)}\left(
\mathrm{sc}_{s}^{0\mathbf{Pth}_{\boldsymbol{\mathcal{B}}^{(2)}}^{\mathbf{f}^{(2)}(1,2)}}\left(
\mathfrak{P}^{(2)}
\right)
\right)
&=
\mathrm{sc}^{([1],2)}_{\boldsymbol{\mathcal{B}}^{(2)}, \varphi(s)}\left(
\mathrm{sc}_{\varphi(s)}^{0\mathbf{Pth}_{\boldsymbol{\mathcal{B}}^{(2)}}}\left(
\mathfrak{P}^{(2)}
\right)
\right)
\tag{1}
\\
&=
\mathrm{sc}_{\varphi(s)}^{0[\mathbf{PT}_{\boldsymbol{\mathcal{B}}^{(1)}}]}\left(
\mathrm{sc}^{([1],2)}_{\boldsymbol{\mathcal{B}}^{(2)}, \varphi(s)}\left(
\mathfrak{P}^{(2)}
\right)
\right)
\tag{2}
\\
&=
\mathrm{sc}_{s}^{0[\mathbf{PT}_{\boldsymbol{\mathcal{B}}^{(1)}}^{\mathbf{f}^{(1)}}]}\left(
\mathrm{sc}^{([1],2)}_{\boldsymbol{\mathcal{B}}^{(2)}, \varphi(s)}\left(
\mathfrak{P}^{(2)}
\right)
\right).
\tag{3}
\end{align*}

The first equality unravels the interpretation of the operation symbol $\mathrm{sc}_{s}^{0}$ in the partial $\Sigma^{\boldsymbol{\mathcal{A}}^{(1)}}$-algebra $\mathbf{Pth}_{\boldsymbol{\mathcal{B}}^{(2)}}^{\mathbf{f}^{(2)}(1,2)}$, introduced in Proposition~\ref{PDPthBCatAlg};
the second equality follows from the fact that, according to Proposition~\ref{PDUCatHom}, $\mathrm{sc}^{([1],2)}_{\boldsymbol{\mathcal{B}}^{(2)}}$ is a $\Lambda^{\boldsymbol{\mathcal{B}}^{(1)}}$-homomorphism from $\mathbf{Pth}_{\boldsymbol{\mathcal{B}}^{(2)}}^{(1,2)}$ to $[\mathbf{PT}_{\boldsymbol{\mathcal{B}}^{(1)}}]$;
finally, the last equality recovers the interpretation of the operation symbol $\mathrm{sc}_{s}^{0}$ in the partial $\Sigma^{\boldsymbol{\mathcal{A}}^{(1)}}$-algebra $[\mathbf{PT}_{\boldsymbol{\mathcal{B}}^{(1)}}^{\mathbf{f}^{(1)}}]$, introduced in Proposition~\ref{PQPTBCatAlg}.

Hence, $\mathrm{sc}^{([1],2)}_{\boldsymbol{\mathcal{B}}^{(2)}, \varphi}$ is compatible with the $0$-source operation.

{\sffamily The mapping $\mathrm{sc}^{([1],2)}_{\boldsymbol{\mathcal{B}}^{(2)}, \varphi}$ is compatible with the $0$-target.}

Let $s$ be a sort in $S$ and let us consider the $0$-target operation symbol $\mathrm{tg}_{s}^{0}$ in $\Sigma^{\boldsymbol{\mathcal{A}}^{(1)}}_{s,s}$. Let $\mathfrak{P}^{(2)}$ be a second-order path in $\mathrm{Pth}_{\boldsymbol{\mathcal{B}}^{(2)}, \varphi(s)}$, then the following equality holds
$$
\mathrm{sc}^{([1],2)}_{\boldsymbol{\mathcal{B}}^{(2)}, \varphi(s)}\left(
\mathrm{tg}_{s}^{0\mathbf{Pth}_{\boldsymbol{\mathcal{B}}^{(2)}}^{\mathbf{f}^{(2)}(1,2)}}\left(
\mathfrak{P}^{(2)}
\right)
\right)
=
\mathrm{tg}_{s}^{0[\mathbf{PT}_{\boldsymbol{\mathcal{B}}^{(1)}}^{\mathbf{f}^{(1)}}]}\left(
\mathrm{sc}^{([1],2)}_{\boldsymbol{\mathcal{B}}^{(2)}, \varphi(s)}\left(
\mathfrak{P}^{(2)}
\right)
\right).
$$

The proof of this case is identical to that of the $0$-source.

Hence, $\mathrm{sc}^{([1],2)}_{\boldsymbol{\mathcal{B}}^{(2)}, \varphi}$ is compatible with the $0$-target operation.

{\sffamily The mapping $\mathrm{sc}^{([1],2)}_{\boldsymbol{\mathcal{B}}^{(2)}, \varphi}$ is compatible with the $0$-composition.}

Let $s$ be a sort in $S$ and let us consider the $0$-composition operation symbol $\circ_{s}^{0}$ in $\Sigma_{ss,s}^{\boldsymbol{\mathcal{A}}^{(1)}}$. Let $\mathfrak{P}^{(2)}$ and $\mathfrak{Q}^{(2)}$ be two second-order paths in $\mathrm{Pth}_{\boldsymbol{\mathcal{B}}^{(2)}, \varphi(s)}$ such that
$$
\mathrm{sc}_{\boldsymbol{\mathcal{B}}^{(2)}, \varphi(s)}^{(0,2)}\left(\mathfrak{Q}^{(2)}\right)
=
\mathrm{tg}_{\boldsymbol{\mathcal{B}}^{(2)}, \varphi(s)}^{(0,2)}\left(\mathfrak{P}^{(2)}\right).
$$

Then the following equality holds
\begin{multline*}
\mathrm{sc}^{([1],2)}_{\boldsymbol{\mathcal{B}}^{(2)}, \varphi(s)}\left(
\mathfrak{Q}^{(2)}
\circ_{s}^{0\mathbf{Pth}_{\boldsymbol{\mathcal{B}}^{(2)}}^{\mathbf{f}^{(2)}(1,2)}}
\mathfrak{P}^{(2)}
\right)
=
\\
\mathrm{sc}^{([1],2)}_{\boldsymbol{\mathcal{B}}^{(2)}, \varphi(s)}\left(
\mathfrak{Q}^{(2)}
\right)
\circ_{s}^{0[\mathbf{PT}_{\boldsymbol{\mathcal{B}}^{(1)}}^{\mathbf{f}^{(1)}}]}
\mathrm{sc}^{([1],2)}_{\boldsymbol{\mathcal{B}}^{(2)}, \varphi(s)}\left(
\mathfrak{P}^{(2)}
\right)
\end{multline*}

The proof of this case is identical to that of the $0$-source.

Hence, $\mathrm{sc}^{([1],2)}_{\boldsymbol{\mathcal{B}}^{(2)}, \varphi}$ is compatible with the $0$-composition operation.

This shows that the mapping $\mathrm{sc}^{([1],2)}_{\boldsymbol{\mathcal{B}}^{(2)}, \varphi}$ is a $\Sigma^{\boldsymbol{\mathcal{A}}^{(1)}}$-homomorphism.

A similar argument applies to the $([1],2)$-target mapping.

This completes the proof.
\end{proof}

\begin{proposition}\label{PIpDUBCatHom}
Let $\mathbf{f}^{(2)}=(\varphi, c, (f^{(i)})_{i\in 3})$ be a second-order morphism from $\boldsymbol{\mathcal{A}}^{(2)}$ to $\boldsymbol{\mathcal{B}}^{(2)}$. Then $\mathrm{ip}^{(2,[1])\sharp}_{\boldsymbol{\mathcal{B}}^{(2)}, \varphi}$ is a $\Sigma^{\boldsymbol{\mathcal{A}}^{(1)}}$-homomorphism from $[\mathbf{PT}_{\boldsymbol{\mathcal{B}}^{(1)}}^{\mathbf{f}^{(1)}}]$ to $\mathbf{Pth}_{\boldsymbol{\mathcal{B}}^{(2)}}^{\mathbf{f}^{(2)}(1,2)}$.
\end{proposition}

\begin{proof}
We prove that $\mathrm{ip}_{\boldsymbol{\mathcal{B}}^{(2)}, \varphi}^{(2,[1])\sharp}$ is compatible with every operation symbol in $\Sigma^{\boldsymbol{\mathcal{A}}^{(1)}}$.

{\sffamily The mapping $\mathrm{ip}_{\boldsymbol{\mathcal{B}}^{(2)}, \varphi}^{(2,[1])\sharp}$ is a $\Sigma$-homomorphism.}

Note that $\mathrm{ip}_{\boldsymbol{\mathcal{B}}^{(2)}, \varphi}^{(2,[1])\sharp} = \mathbf{c}_{\mathfrak{d}}^{\ast} (\mathrm{ip}_{\boldsymbol{\mathcal{B}}^{(2)}}^{(2,[1])\sharp})$. By Proposition~\ref{PDUIpCatHom}, the mapping $\mathrm{ip}_{\boldsymbol{\mathcal{B}}^{(2)}}^{(2,[1])\sharp}$ is a $\Lambda^{\boldsymbol{\mathcal{B}}^{(1)}}$-homomorphism, thus in particular a $\Lambda$-homomorphism. Therefore, it follows from Proposition~\ref{PFunSig} that the mapping $\mathrm{ip}_{\boldsymbol{\mathcal{B}}^{(2)}, \varphi}^{(2,[1])\sharp}$ is a $\Sigma$-homomorphism.

{\sffamily The mapping $\mathrm{ip}_{\boldsymbol{\mathcal{B}}^{(2)}, \varphi}^{(2,[1])\sharp}$ is compatible with the first-order rewrite rules.}

Let $s$ be a sort in $S$ and $\mathfrak{p}$ a rewrite rule in $\mathcal{A}_{s}^{(1)}$. Thus, the following chain of equalities holds
\begin{flushleft}
$
\mathrm{ip}_{\boldsymbol{\mathcal{B}}^{(2)}, \varphi(s)}^{(2,[1])\sharp}\left(
\mathfrak{p}^{[\mathbf{PT}_{\boldsymbol{\mathcal{B}}^{(1)}}^{\mathbf{f}^{(1)}}]}
\right)
$
\allowdisplaybreaks
\begin{align*}
&=
\mathrm{ip}_{\boldsymbol{\mathcal{B}}^{(2)}, \varphi(s)}^{(2,[1])\sharp}\left(
\left[
\mathrm{CH}^{(1)}_{\boldsymbol{\mathcal{B}}^{(1)}, \varphi(s)} \left(
f^{(1)\flat}_{s}\left(
\mathfrak{p}^{\mathbf{Pth}_{\boldsymbol{\mathcal{A}}^{(1)}}}
\right)
\right)
\right]_{\varphi(s)}
\right)
\tag{1}
\\
&=
\mathrm{ip}_{\boldsymbol{\mathcal{B}}^{(2)}, \varphi(s)}^{(2,[1])\sharp}\left(
\left[
\mathrm{CH}^{(1)}_{\boldsymbol{\mathcal{B}}^{(1)}, \varphi(s)} \left(
f^{(1)\flat}_{s}\left(
\mathrm{ech}^{(1,\mathcal{A}^{(1)})}_{\boldsymbol{\mathcal{A}}^{(1)}, s}\left(
\mathfrak{p}
\right)
\right)
\right)
\right]_{\varphi(s)}
\right)
\tag{2}
\\
&=
\mathrm{ip}_{\boldsymbol{\mathcal{B}}^{(2)}, \varphi(s)}^{(2,[1])\sharp}\left(
\left[
\mathrm{CH}^{(1)}_{\boldsymbol{\mathcal{B}}^{(1)}, \varphi(s)} \left(
f^{(1)\flat}_{s}\left(
\mathrm{ip}^{(1,X)@}_{\boldsymbol{\mathcal{A}}^{(1)}, s}\left(
\eta^{(1,\mathcal{A}^{(1)})}_{\boldsymbol{\mathcal{A}}^{(1)}, s}\left(
\mathfrak{p}
\right)
\right)
\right)
\right)
\right]_{\varphi(s)}
\right)
\tag{3}
\\
&=
\mathrm{ip}_{\boldsymbol{\mathcal{B}}^{(2)}, \varphi(s)}^{(2,[1])\sharp}\left(
f^{[1]@}_{s}\left(
\left[
\eta^{(1,\mathcal{A}^{(1)})}_{\boldsymbol{\mathcal{A}}^{(1)}, s}\left(
\mathfrak{p}
\right)
\right]_{s}
\right)
\right)
\tag{4}
\\
&=
\mathrm{ip}_{\boldsymbol{\mathcal{B}}^{(2)}, \varphi(s)}^{(2,[1])\sharp}\left(
f^{[1]@}_{s}\left(
\eta^{([1],\mathcal{A}^{(1)})}_{\boldsymbol{\mathcal{A}}^{(1)}, s}\left(
\mathfrak{p}
\right)
\right)
\right)
\tag{5}
\\
&=
f^{2\flat}_{s}\left(
\mathrm{ip}_{\boldsymbol{\mathcal{A}}^{(2)}, s}^{(2,[1])\sharp}\left(
\eta^{([1],\mathcal{A}^{(1)})}_{\boldsymbol{\mathcal{A}}^{(1)}, s}\left(
\mathfrak{p}
\right)
\right)
\right)
\tag{6}
\\
&=
f^{2\flat}_{s}\left(
\mathrm{ech}^{(2,\mathcal{A}^{(1)})}_{\boldsymbol{\mathcal{A}}^{(2)}, s}\left(
\mathfrak{p}
\right)
\right)
\tag{7}
\\
&=
f^{2\flat}_{s}\left(
\mathfrak{p}^{\mathbf{Pth}_{\boldsymbol{\mathcal{A}}^{(2)}}}
\right)
\tag{8}
\\
&=
\mathfrak{p}^{\mathbf{Pth}_{\boldsymbol{\mathcal{B}}^{(2)}}^{\mathbf{f}^{(2)}(1,2)}}.
\tag{9}
\end{align*}
\end{flushleft}

In the just stated chain of equalities, the first equality unravels the interpretation of the constant symbol $\mathfrak{p}$ in the partial $\Sigma^{\boldsymbol{\mathcal{A}}^{(1)}}$-algebra $[\mathbf{PT}_{\boldsymbol{\mathcal{B}}^{(1)}}^{\mathbf{f}^{(1)}}]$, introduced in Proposition~\ref{PPTBCatAlg};
the second equality unravels the interpretation of the constant symbol $\mathfrak{p}$ in the partial $\Sigma^{\boldsymbol{\mathcal{A}}^{(1)}}$-algebra $\mathbf{Pth}_{\boldsymbol{\mathcal{A}}^{(1)}}$, introduced in Proposition~\ref{PPthCatAlg};
the third equality follows from the fact that, according to Proposition~\ref{PIpEch}, the mapping $\mathrm{ech}^{(1,\mathcal{A}^{(1)})}_{\boldsymbol{\mathcal{A}}^{(1)}}$ is equal to $\mathrm{ip}^{(1,X)@}_{\boldsymbol{\mathcal{A}}^{(1)}} \circ \eta^{1, \mathcal{A}^{(1)})}_{\boldsymbol{\mathcal{A}}^{(1)}}$;
the fourth equality recovers the definition of the mapping $f^{[1]@}$ at a path term class, introduced in Definition~\ref{DQPthExt};
the fifth equality recovers the definition of the mapping $\eta^{([1], \mathcal{A}^{(1)})}_{\boldsymbol{\mathcal{A}}^{(1)}}$, introduced in Definition~\ref{DPTQEta};
the sixth equality follows from Proposition~\ref{PDPthExt};
the seventh equality follows from Proposition~\ref{PDBasicEqEch};
the eighth equality recovers the interpretation of the constant symbol $\mathfrak{p}$ in the partial $\Sigma^{\boldsymbol{\mathcal{A}}^{(1)}}$-algebra $\mathbf{Pth}_{\boldsymbol{\mathcal{A}}^{(2)}}$, introduced in Proposition~\ref{PDPthCatAlg};
finally, the last equality recovers the interpretation of the constant symbol $\mathfrak{p}$ in the partial $\Sigma^{\boldsymbol{\mathcal{A}}^{(1)}}$-algebra $\mathbf{Pth}_{\boldsymbol{\mathcal{B}}^{(2)}}^{\mathbf{f}^{(2)}(1,2)}$, introduced in Proposition~\ref{PDPthBCatAlg};

Hence, $\mathrm{ip}_{\boldsymbol{\mathcal{B}}^{(2)}, \varphi}^{(2,[1])\sharp}$ is compatible with the first-order rewrite rules.

{\sffamily The mapping $\mathrm{ip}_{\boldsymbol{\mathcal{B}}^{(2)}, \varphi}^{(2,[1])\sharp}$ is compatible with the $0$-source.}

Let $s$ be a sort in $S$ and let us consider the $0$-source operation symbol $\mathrm{sc}_{s}^{0}$ in $\Sigma^{\boldsymbol{\mathcal{A}}^{(1)}}_{s,s}$. Let $\mathfrak{P}^{(2)}$ be a second-order path in $\mathrm{Pth}_{\boldsymbol{\mathcal{B}}^{(2)}, \varphi(s)}$.

The following chain of equalities holds
\allowdisplaybreaks
\begin{align*}
\mathrm{ip}_{\boldsymbol{\mathcal{B}}^{(2)}, \varphi(s)}^{(2,[1])\sharp}\left(
\mathrm{sc}_{s}^{0[\mathbf{PT}_{\boldsymbol{\mathcal{B}}^{(1)}}^{\mathbf{f}^{(1)}}]}\left(
\mathfrak{P}^{(2)}
\right)
\right)
&=
\mathrm{ip}_{\boldsymbol{\mathcal{B}}^{(2)}, \varphi(s)}^{(2,[1])\sharp}\left(
\mathrm{sc}_{\varphi(s)}^{0[\mathbf{PT}_{\boldsymbol{\mathcal{B}}^{(1)}}]}\left(
\mathfrak{P}^{(2)}
\right)
\right)
\tag{1}
\\
&=
\mathrm{sc}_{\varphi(s)}^{0\mathbf{Pth}_{\boldsymbol{\mathcal{B}}^{(2)}}}\left(
\mathrm{ip}_{\boldsymbol{\mathcal{B}}^{(2)}, \varphi(s)}^{(2,[1])\sharp}\left(
\mathfrak{P}^{(2)}
\right)
\right)
\tag{2}
\\
&=
\mathrm{sc}_{s}^{0\mathbf{Pth}_{\boldsymbol{\mathcal{B}}^{(2)}}^{\mathbf{f}^{(2)}}}\left(
\mathrm{ip}_{\boldsymbol{\mathcal{B}}^{(2)}, \varphi(s)}^{(2,[1])\sharp}\left(
\mathfrak{P}^{(2)}
\right)
\right).
\tag{3}
\end{align*}

The first equality unravels the interpretation of the operation symbol $\mathrm{sc}_{s}^{0}$ in the partial $\Sigma^{\boldsymbol{\mathcal{A}}^{(1)}}$-algebra $[\mathbf{PT}_{\boldsymbol{\mathcal{B}}^{(1)}}^{\mathbf{f}^{(1)}}]$, introduced in Proposition~\ref{PQPTBCatAlg};
the second equality follows from the fact that, according to Proposition~\ref{PDUIpCatHom}, $\mathrm{ip}^{(2,[1])\sharp}_{\boldsymbol{\mathcal{B}}^{(2)}}$ is a $\Lambda^{\boldsymbol{\mathcal{B}}^{(1)}}$-homomorphism from $[\mathbf{PT}_{\boldsymbol{\mathcal{B}}^{(2)}}^{\mathbf{f}^{(2)}}]$ to $\mathbf{Pth}_{\boldsymbol{\mathcal{B}}^{(2)}}^{(1,2)}$;
finally, the last equality recovers the interpretation of the operation symbol $\mathrm{sc}_{s}^{0}$ in the partial $\Sigma^{\boldsymbol{\mathcal{A}}^{(1)}}$-algebra $\llbracket\mathbf{Pth}_{\boldsymbol{\mathcal{B}}^{(2)}}^{\mathbf{f}^{(2)}(1,2)}\rrbracket$, introduced in Proposition~\ref{PDQPthBCatAlg}.

Hence, $\mathrm{ip}_{\boldsymbol{\mathcal{B}}^{(2)}, \varphi}^{(2,[1])\sharp}$ is compatible with the $0$-source operation.

{\sffamily The mapping $\mathrm{ip}_{\boldsymbol{\mathcal{B}}^{(2)}, \varphi}^{(2,[1])\sharp}$ is compatible with the $0$-target.}

Let $s$ be a sort in $S$ and let us consider the $0$-target operation symbol $\mathrm{tg}_{s}^{0}$ in $\Sigma^{\boldsymbol{\mathcal{A}}^{(1)}}_{s,s}$. Let $\mathfrak{P}^{(2)}$ be a second-order path in $\mathrm{Pth}_{\boldsymbol{\mathcal{B}}^{(2)}, \varphi(s)}$, then the following equality holds
$$
\mathrm{ip}_{\boldsymbol{\mathcal{B}}^{(2)}, \varphi(s)}^{(2,[1])\sharp}\left(
\mathrm{tg}_{s}^{0[\mathbf{PT}_{\boldsymbol{\mathcal{B}}^{(1)}}^{\mathbf{f}^{(1)}}]}\left(
\mathfrak{P}^{(2)}
\right)
\right)
=
\mathrm{tg}_{s}^{0\mathbf{Pth}_{\boldsymbol{\mathcal{B}}^{(2)}}^{\mathbf{f}^{(2)}(1,2)}}\left(
\mathrm{ip}_{\boldsymbol{\mathcal{B}}^{(2)}, \varphi(s)}^{(2,[1])\sharp}\left(
\mathfrak{P}^{(2)}
\right)
\right).
$$

The proof of this case is identical to that of the $0$-source.

Hence, $\mathrm{ip}_{\boldsymbol{\mathcal{B}}^{(2)}, \varphi}^{(2,[1])\sharp}$ is compatible with the $0$-target operation.

{\sffamily The mapping $\mathrm{ip}_{\boldsymbol{\mathcal{B}}^{(2)}, \varphi}^{(2,[1])\sharp}$ is compatible with the $0$-composition.}

Let $s$ be a sort in $S$ and let us consider the $0$-composition operation symbol $\circ_{s}^{0}$ in $\Sigma_{ss,s}^{\boldsymbol{\mathcal{A}}^{(1)}}$. Let $\mathfrak{P}^{(2)}$ and $\mathfrak{Q}^{(2)}$ be two second-order paths in $\mathrm{Pth}_{\boldsymbol{\mathcal{B}}^{(2)}, \varphi(s)}$ such that
$$
\mathrm{sc}_{\boldsymbol{\mathcal{B}}^{(2)}, \varphi(s)}^{(0,2)}\left(\mathfrak{Q}^{(2)}\right)
=
\mathrm{tg}_{\boldsymbol{\mathcal{B}}^{(2)}, \varphi(s)}^{(0,2)}\left(\mathfrak{P}^{(2)}\right).
$$

Then the following equality holds
\begin{multline*}
\mathrm{ip}_{\boldsymbol{\mathcal{B}}^{(2)}, \varphi(s)}^{(2,[1])\sharp}\left(
\mathfrak{Q}^{(2)}
\circ_{s}^{0[\mathbf{PT}_{\boldsymbol{\mathcal{B}}^{(1)}}^{\mathbf{f}^{(1)}}]}
\mathfrak{P}^{(2)}
\right)
=
\\
\mathrm{ip}_{\boldsymbol{\mathcal{B}}^{(2)}, \varphi(s)}^{(2,[1])\sharp}\left(
\mathfrak{Q}^{(2)}
\right)
\circ_{s}^{0\mathbf{Pth}_{\boldsymbol{\mathcal{B}}^{(2)}}^{\mathbf{f}^{(2)}(1,2)}}
\mathrm{ip}_{\boldsymbol{\mathcal{B}}^{(2)}, \varphi(s)}^{(2,[1])\sharp}\left(
\mathfrak{P}^{(2)}
\right)
\end{multline*}

The proof of this case is identical to that of the $0$-source.

Hence, $\mathrm{ip}_{\boldsymbol{\mathcal{B}}^{(2)}, \varphi}^{(2,[1])\sharp}$ is compatible with the $0$-composition operation.

This shows that the mapping $\mathrm{ip}_{\boldsymbol{\mathcal{B}}^{(2)}, \varphi}^{(2,[1])\sharp}$ is a $\Sigma^{\boldsymbol{\mathcal{A}}^{(1)}}$-homomorphism.

This completes the proof.
\end{proof}

\begin{proposition}
\label{PDPthExtCatHom}
Let $\mathbf{f}^{(2)}=(\varphi, c, (f^{(i)})_{i\in 3})$ be a second-order morphism from $\boldsymbol{\mathcal{A}}^{(2)}$ to $\boldsymbol{\mathcal{B}}^{(2)}$. The second-order path extension mapping $f^{(2)\flat}$ is a $\Sigma^{\boldsymbol{\mathcal{A}}^{(1)}}$-homomorphism from $\mathbf{Pth}_{\boldsymbol{\mathcal{A}}^{(2)}}^{(1,2)}$ to $\mathbf{Pth}_{\boldsymbol{\mathcal{B}}^{(2)}}^{\mathbf{f}^{(2)}(1,2)}$. 
\end{proposition}

\begin{proof}
We prove that $f^{(2)\flat}$ is compatible with every operation symbol in $\Sigma^{\boldsymbol{\mathcal{A}}^{(1)}}$.

{\sffamily The mapping $f^{(2)\flat}$ is a $\Sigma$-homomorphism.}

It follows from Proposition~\ref{PDPthExtHom}, that $f^{(2)\flat}$ is a $\Sigma$-homomorphism .

{\sffamily The mapping $f^{(2)\flat}$ is compatible with the first-order rewrite rules.}

For every sort $s$ in $S$ and first-order rewrite rule $\mathfrak{p}$ in $\mathcal{A}^{(1)}_{s}$, the following chain of equaities holds
\begin{align*}
f^{(2)\flat}_{s}\left(
\mathfrak{p}^{\mathbf{Pth}_{\boldsymbol{\mathcal{A}}^{(2)}}^{(1,2)}}
\right)
&=
f^{(2)\flat}_{s}\left(
\mathrm{ech}^{(2,\mathcal{A}^{(1)})}_{\boldsymbol{\mathcal{A}}^{(2)}, s}\left(
\mathfrak{p}
\right)
\right)
\tag{1}
\\
&=
f^{(2)\flat}_{s}\left(
\mathrm{ip}^{(2,[1])\sharp}_{\boldsymbol{\mathcal{A}}^{(2)}, s} \left(
\eta^{([1],\mathcal{A}^{(1)})}_{\boldsymbol{\mathcal{A}}^{(1)}, s}\left(
\mathfrak{p}
\right)
\right)
\right)
\tag{2}
\\
&=
\mathrm{ip}^{(2,[1])\sharp}_{\boldsymbol{\mathcal{B}}^{(2)}, \varphi(s)} \left(
f^{[1]@}_{s}\left(
\eta^{([1],\mathcal{A}^{(1)})}_{\boldsymbol{\mathcal{A}}^{(1)}, s}\left(
\mathfrak{p}
\right)
\right)
\right)
\tag{3}
\\
&=
\mathrm{ip}^{(2,[1])\sharp}_{\boldsymbol{\mathcal{B}}^{(2)}, \varphi(s)} \left(
f^{[1]@}_{s}\left(
\mathfrak{p}^{[\mathbf{PT}_{\boldsymbol{\mathcal{A}}^{(1)}}]}
\right)
\right)
\tag{4}
\\
&=
\mathrm{ip}^{(2,[1])\sharp}_{\boldsymbol{\mathcal{B}}^{(2)}, \varphi(s)} \left(
\mathfrak{p}^{[\mathbf{PT}_{\boldsymbol{\mathcal{B}}^{(1)}}^{\mathbf{f}^{(1)}}]}
\right)
\tag{5}
\\
&=
\mathfrak{p}^{\mathbf{Pth}_{\boldsymbol{\mathcal{B}}^{(2)}}^{\mathbf{f}^{(2)}(1,2)}}
\tag{6}
\end{align*}

The first equality unravels the definition of the constant operation symbol $\mathfrak{p}$ in the partial $\Sigma^{\boldsymbol{\mathcal{A}}^{(1)}}$-algebra $\mathbf{Pth}_{\boldsymbol{\mathcal{A}}^{(2)}}^{(1,2)}$, introduced in Proposition~\ref{PDPthCatAlg};
the second equality follows from Proposition~\ref{PDBasicEqEch};
the third equality follows from Proposition~\ref{PDPthExt};
the fourth equality recovers the definition of the constant operation symbol $\mathfrak{p}$ in the partial $\Sigma^{\boldsymbol{\mathcal{A}}^{(1)}}$-algebra $[\mathbf{PT}_{\boldsymbol{\mathcal{A}}^{(1)}}]$, introduced in Proposition~\ref{PPTQCatAlg};
the fifth equality follows from the fact that, according to Definition~\ref{DQPthExt}, the quotient path-extension mapping is a $\Sigma^{\boldsymbol{\mathcal{A}}^{(1)}}$-homomorphism;
finally, the last equality follows from the fact that, according to Proposition~\ref{PIpDUBCatHom}, $\mathrm{ip}^{(2,[1])\sharp}_{\boldsymbol{\mathcal{B}}^{(2)}, \varphi}$ is a $\Sigma^{\boldsymbol{\mathcal{A}}^{(1)}}$-homomorphism.

Hence, $f^{(2)\flat}$ is compatible with the first-order rewrite rules.

{\sffamily The mapping $f^{(2)\flat}$ is compatible with the $0$-source.}

For every sort $s$ in $S$, the following chain of equalities holds
\allowdisplaybreaks
\begin{align*}
f^{(2)\flat}_{s} \circ \mathrm{sc}^{0\mathbf{Pth}_{\boldsymbol{\mathcal{A}}^{(2)}}^{(1,2)}}_{s}
&=
f^{(2)\flat}_{s} \circ \mathrm{ip}^{(2,0)\sharp}_{\boldsymbol{\mathcal{A}}^{(2)}, s} \circ \mathrm{sc}^{(0,2)}_{\boldsymbol{\mathcal{A}}^{(2)}, s}
\tag{1}
\\
&=
\mathrm{ip}^{(2,0)\sharp}_{\boldsymbol{\mathcal{B}}^{(2)}, \varphi(s)} \circ f^{(0)\sharp}_{s} \circ \mathrm{sc}^{(0,2)}_{\boldsymbol{\mathcal{A}}^{(2)}, s}
\tag{2}
\\
&=
\mathrm{ip}^{(2,0)\sharp}_{\boldsymbol{\mathcal{B}}^{(2)}, \varphi(s)} \circ \mathrm{sc}^{(0,2)}_{\boldsymbol{\mathcal{B}}^{(2)}, \varphi(s)} \circ f^{(2)\flat}_{s}
\tag{3}
\\
&=
\mathrm{sc}^{0\mathbf{Pth}_{\boldsymbol{\mathcal{B}}^{(2)}}^{(1,2)}}_{\varphi(s)} \circ f^{(2)\flat}_{\varphi(s)}
\tag{4}
\\
&=
\mathrm{sc}^{0\mathbf{Pth}_{\boldsymbol{\mathcal{B}}^{(2)}}^{\mathbf{f}^{(2)}(1,2)}}_{s} \circ f^{(2)\flat}_{s}.
\tag{5}
\end{align*}

The first equality unravels the definition of the interpretation of the operation symbol $\mathrm{sc}_{s}^{0}$ in the partial $\Sigma^{\boldsymbol{\mathcal{A}}^{(1)}}$-algebra $\mathbf{Pth}_{\boldsymbol{\mathcal{A}}^{(2)}}^{(1,2)}$, introduced in Proposition~\ref{PDPthCatAlg};
the second equality follows from Proposition~\ref{PDPthExtDZIp};
the third equality follows from Proposition~\ref{PDPthExtDZScTg};
the fourth equality recovers the definition of the interpretation of the operation symbol $\mathrm{sc}_{s}^{0}$ in the partial $\Lambda^{\boldsymbol{\mathcal{B}}^{(1)}}$-algebra $\mathbf{Pth}_{\boldsymbol{\mathcal{B}}^{(2)}}^{(1,2)}$, introduced in Proposition~\ref{PDPthCatAlg};
finally, the last equality recovers the definition of the interpretation of the operation symbol $\mathrm{sc}_{s}^{0}$ in the partial $\Sigma^{\boldsymbol{\mathcal{A}}^{(1)}}$-algebra $\mathbf{Pth}_{\boldsymbol{\mathcal{B}}^{(2)}}^{\mathbf{f}^{(2)}(1,2)}$, introduced in Proposition~\ref{PDPthBCatAlg};

Hence, $f^{(2)\flat}$ is compatible with the $0$-source operation.

{\sffamily The mapping $f^{(2)\flat}$ is compatible with the $0$-target.}

The proof is similar to that of the $0$-source operation.

Hence, $f^{(2)\flat}$ is compatible with the $0$-target operation.

{\sffamily The mapping $f^{(2)\flat}$ is compatible with the $0$-composition.}

Let $s$ be a sort in $S$ and let $\mathfrak{P}^{(2)}_{1}$ and $\mathfrak{P}^{(2)}_{0}$ second-order paths in $\mathrm{Pth}_{\boldsymbol{\mathcal{A}}^{(2)}}$ such that their $0$-composition
$$
\mathfrak{P}^{(2)}
=
\mathfrak{P}^{(2)}_{1} \circ_{s}^{0\mathbf{Pth}_{\boldsymbol{\mathcal{A}}^{(1)}}^{(1,2)}} \mathfrak{P}^{(2)}_{0}
$$
is defined. Then note that the path $\mathfrak{P}^{(2)}$ is a coherent head-constant echelonless second-order path associated with the operation symbol $\circ_{s}^{0}$. Thus, according to the definition of the second-order path extension mapping, introduced in Proposition~\ref{PDPthExt},
\begin{align*}
f^{(2)\flat}_{s}\left(
\mathfrak{P}^{(2)}
\right)
&=
f^{(2)\flat}_{s}\left(
\mathfrak{P}^{(2)}_{1}
\right)
\circ_{\varphi(s)}^{0\mathbf{Pth}_{\boldsymbol{\mathcal{B}}^{(2)}}}
f^{(2)\flat}_{s}\left(
\mathfrak{P}^{(2)}_{0}
\right)
\tag{1}
\\
&=
f^{(2)\flat}_{s}\left(
\mathfrak{P}^{(2)}_{1}
\right)
\circ_{s}^{0\mathbf{Pth}_{\boldsymbol{\mathcal{B}}^{(2)}}^{\mathbf{f}^{(2)}(1,2)}}
f^{(2)\flat}_{s}\left(
\mathfrak{P}^{(2)}_{0}
\right).
\tag{2}
\end{align*}

Hence, $f^{(2)\flat}$ is compatible with the $0$-composition operation.
\end{proof}
\chapter{Derived partial $\Sigma^{\boldsymbol{\mathcal{A}}^{(2)}}$-algebras}\label{S3F}

In this chapter we define a structure of partial $\Sigma^{\boldsymbol{\mathcal{A}}^{(2)}}$-algebra on $\mathrm{Pth}_{\boldsymbol{\mathcal{B}}^{(2)}, \varphi}$ and $\llbracket \mathrm{Pth}_{\boldsymbol{\mathcal{B}}^{(2)}} \rrbracket_{\varphi}$, which we will denote by $\mathbf{Pth}_{\boldsymbol{\mathcal{B}}^{(2)}}^{\mathbf{f}^{(2)}}$ and $\llbracket \mathbf{Pth}_{\boldsymbol{\mathcal{B}}^{(2)}}^{\mathbf{f}^{(2)}} \rrbracket$, respectively. We show that the second-order path extension mapping is not a $\Sigma^{\boldsymbol{\mathcal{A}}^{(2)}}$-homomorphism. However, we show that the projection mapping $\mathrm{pr}^{\llbracket\cdot\rrbracket}_{\boldsymbol{\mathcal{B}}^{(2)}, \varphi}$ is a $\Sigma^{\boldsymbol{\mathcal{A}}^{(2)}}$-homomorphism. We next define a structure of partial $\Sigma^{\boldsymbol{\mathcal{A}}^{(2)}}$-algebra on $\mathrm{PT}_{\boldsymbol{\mathcal{B}}^{(2)}, \varphi}$ and $\llbracket\mathrm{PT}_{\boldsymbol{\mathcal{B}}^{(2)}}\rrbracket_{\varphi}$, which we will denote by $\mathbf{PT}_{\boldsymbol{\mathcal{B}}^{(2)}}^{\mathbf{f}^{(2)}}$ and $\llbracket \mathbf{PT}_{\boldsymbol{\mathcal{B}}^{(2)}}^{\mathbf{f}^{(2)}} \rrbracket$, respectively. We then show that the projection mapping $\mathrm{pr}^{\Theta^{\llbracket2\rrbracket}}_{\boldsymbol{\mathcal{B}}^{(2)}, \varphi}$ is a $\Sigma^{\boldsymbol{\mathcal{A}}^{(2)}}$-homomorphism. Moreover, we show that $\llbracket \mathbf{Pth}_{\boldsymbol{\mathcal{B}}^{(2)}}^{\mathbf{f}^{(2)}} \rrbracket$ and $\llbracket \mathbf{PT}_{\boldsymbol{\mathcal{B}}^{(2)}}^{\mathbf{f}^{(2)}} \rrbracket$ are isomorphic partial $\Sigma^{\boldsymbol{\mathcal{A}}^{(2)}}$-algebras. Finally, we define a structure of partial $\Sigma^{\boldsymbol{\mathcal{A}}^{(2)}}$-algebra on $\T_{\boldsymbol{\mathcal{E}}^{\boldsymbol{\mathcal{B}}^{(2)}}}(\mathbf{Pth}_{\boldsymbol{\mathcal{B}}^{(2)}})_{\varphi}$, which we will denote by $\mathbf{T}_{\boldsymbol{\mathcal{E}}^{\boldsymbol{\mathcal{B}}^{(2)}}}^{\mathbf{f}^{(2)}}(\mathbf{Pth}_{\boldsymbol{\mathcal{B}}^{(2)}})$ and we prove that the partial $\Sigma^{\boldsymbol{\mathcal{A}}^{(2)}}$-algebras $\llbracket \mathbf{Pth}_{\boldsymbol{\mathcal{B}}^{(2)}}^{\mathbf{f}^{(2)}} \rrbracket$ and $\mathbf{T}_{\boldsymbol{\mathcal{E}}^{\boldsymbol{\mathcal{B}}^{(2)}}}^{\mathbf{f}^{(2)}}(\mathbf{Pth}_{\boldsymbol{\mathcal{B}}^{(2)}})$ are also isomorphic.

\section{A structure of partial $\Sigma^{\boldsymbol{\mathcal{A}}^{(2)}}$-algebra on $\mathrm{Pth}_{\boldsymbol{\mathcal{B}}^{(2)}, \varphi}$}

We next show that, given a second-order morphism $\mathbf{f}^{(2)}$, $\mathrm{Pth}_{\boldsymbol{\mathcal{B}}^{(2)}, \varphi}$ is equipped in a natural way with a structure of partial $\Sigma^{\boldsymbol{\mathcal{A}}^{(2)}}$-algebra that we denote by $\mathbf{Pth}_{\boldsymbol{\mathcal{B}}^{(2)}}^{\mathbf{f}^{(2)}}$.

\begin{proposition}
\label{PDPthBDCatAlg}
Let $\mathbf{f}^{(2)}=(\varphi, c, (f^{(i)})_{i\in 3})$ be a second-order morphism from $\boldsymbol{\mathcal{A}}^{(2)}$ to $\boldsymbol{\mathcal{B}}^{(2)}$. Then the $S$-sorted set $\mathrm{Pth}_{\boldsymbol{\mathcal{B}}^{(2)}, \varphi}$ is equipped, in a natural way, with a structure of partial $\Sigma^{\boldsymbol{\mathcal{A}}^{(2)}}$-algebra.
\end{proposition}

\begin{proof}
Let us denote by $\mathbf{Pth}_{\boldsymbol{\mathcal{B}}^{(2)}}^{\mathbf{f}^{(2)}}$ the $\Sigma^{\boldsymbol{\mathcal{A}}^{(2)}}$-algebra defined as follows:

\textsf{(1)}
The underlying $S$-sorted set of $\mathbf{Pth}_{\boldsymbol{\mathcal{B}}^{(2)}}^{\mathbf{f}^{(2)}}$ is $\mathrm{Pth}_{\boldsymbol{\mathcal{B}}^{(2)}, \varphi}=(\mathrm{Pth}_{\boldsymbol{\mathcal{B}}^{(2)}, \varphi(s)})_{s\in S}$.

\textsf{(2)}
For every $(\mathbf{s}, s)\in S^{\ast}\times S$ and every operation symbol $\sigma\in\Sigma_{\mathbf{s},s}$, the operation $\sigma^{\mathbf{Pth}_{\boldsymbol{\mathcal{B}}^{(2)}}^{\mathbf{f}^{(2)}}}$ is given by the interpretation of $\sigma$ in the $\Sigma$-algebra $\mathbf{Pth}_{\boldsymbol{\mathcal{B}}^{(2)}}^{\mathbf{f}^{(2)}(0,2)}$ introduced in Remark~\ref{RDDSigmaAlg}.

\textsf{(3)}
For every $s\in S$ and every $\mathfrak{p}\in\mathcal{A}_{s}^{(1)}$, the constant $\mathfrak{p}^{\mathbf{Pth}_{\boldsymbol{\mathcal{B}}^{(2)}}^{\mathbf{f}^{(2)}}}$ is given by the interpretation of $\mathrm{p}$ in $\mathbf{Pth}_{\boldsymbol{\mathcal{A}}^{(2)}}^{\mathbf{f}^{(2)}(1,2)}$ that , we recall, was stated in Proposition~\ref{PDPthBCatAlg}.

\textsf{(4)}
For every $s\in S$, the interpretations of the operations $\mathrm{sc}_{s}^{0}$ and $\mathrm{tg}_{s}^{0}$ are given by the interpretation of $\mathrm{sc}_{s}^{0}$ and $\mathrm{tg}_{s}^{0}$ in $\mathbf{Pth}_{\boldsymbol{\mathcal{A}}^{(2)}}^{\mathbf{f}^{(2)}(1,2)}$ that , we recall, was stated in Proposition~\ref{PDPthBCatAlg}.

\textsf{(5)}
Similarly,  for every $s \in S$, the interpretations of the operation $\circ_{s}^{0}$ is given by the interpretation of $\circ_{s}^{0}$ in $\mathbf{Pth}_{\boldsymbol{\mathcal{A}}^{(2)}}^{\mathbf{f}^{(2)}(1,2)}$ that , we recall, was stated in Proposition~\ref{PDPthBCatAlg}.

\textsf{(6)}
For every $s\in S$ and every $\mathfrak{p}^{(2)}\in\mathcal{A}_{s}^{(2)}$, the constant $\mathfrak{p}^{(2)\mathbf{Pth}_{\boldsymbol{\mathcal{B}}^{(2)}}^{\mathbf{f}^{(2)}}}$ is given by the image of $\mathfrak{p}^{(2)}$ under the mapping $f^{(2)}$, i.e., we set

$$
\mathfrak{p}^{(2)\mathbf{Pth}_{\boldsymbol{\mathcal{B}}^{(2)}}^{\mathbf{f}^{(2)}}}
=
f^{(2)}_{s}\left(
\mathfrak{p}^{(2)}
\right)
\in
\mathrm{Pth}_{\boldsymbol{\mathcal{B}}^{(2)}, \varphi(s)}
$$.

\textsf{(7)}
For every $s\in S$, the interpretations of the operations $\mathrm{sc}_{s}^{1}$ and $\mathrm{tg}_{s}^{1}$ are given by
\allowdisplaybreaks
\begin{align*}
\mathrm{sc}_{s}^{1\mathbf{Pth}_{\boldsymbol{\mathcal{B}}^{(2)}}^{\mathbf{f}^{(2)}}}
&=
\mathrm{sc}_{\varphi(s)}^{1\mathbf{Pth}_{\boldsymbol{\mathcal{B}}^{(2)}}}
&&\mbox{and}&
\mathrm{tg}_{s}^{1\mathbf{Pth}_{\boldsymbol{\mathcal{B}}^{(2)}}^{\mathbf{f}^{(2)}}}
&=
\mathrm{tg}_{\varphi(s)}^{1\mathbf{Pth}_{\boldsymbol{\mathcal{B}}^{(2)}}}.
\end{align*}
That is, their interpretations are given by the interpretations of $\mathrm{sc}_{\varphi(s)}^{1}$ and $\mathrm{tg}_{\varphi(s)}^{1}$ in $\mathbf{Pth}_{\boldsymbol{\mathcal{B}}^{(2)}}$ that, we recall, was introduced in Proposition~\ref{PDPthDCatAlg}.

\textsf{(8)}
Similarly, for every $s \in S$,  the interpretation of the partial binary operation $\circ_{s}^{1}$ is given by
$$
\circ_{s}^{1\mathbf{Pth}_{\boldsymbol{\mathcal{B}}^{(2)}}^{\mathbf{f}^{(2)}}}
=
\circ_{\varphi(s)}^{1\mathbf{Pth}_{\boldsymbol{\mathcal{B}}^{(2)}}}.
$$
That is, its interpretation is given by the interpretation of $\circ_{\varphi(s)}^{1}$ in $\mathbf{Pth}_{\boldsymbol{\mathcal{B}}^{(2)}}$ that, we recall, was introduced in Proposition~\ref{PDPthDCatAlg}.

This completes the definition of the partial $\Sigma^{\boldsymbol{\mathcal{A}}^{(2)}}$-algebra $\mathbf{Pth}_{\boldsymbol{\mathcal{B}}^{(2)}}^{\mathbf{f}^{(2)}}$.
\end{proof}

Despite the relations of the second-order path extension mapping with the source, target, and identity path mappings, it cannot be inferred that it is a $\Sigma^{\boldsymbol{\mathcal{A}}^{(2)}}$-homomorphism. This is largely due to the relationship between the $1$-composition of paths and the Artinian order used to define the mapping $f^{(2)\flat}$. Below we present a counterexample to the assertion that $f^{(2)\flat}$ is a $\Sigma^{\boldsymbol{\mathcal{A}}^{(2)}}$-homomorphism.

\begin{remark}
Let $\mathbf{f}^{(2)}=(\varphi, c, (f^{(i)})_{i\in 3})$ be a second-order morphism from $\boldsymbol{\mathcal{A}}^{(2)}$ to $\boldsymbol{\mathcal{B}}^{(2)}$. The second-order path extension mapping $f^{(2)\flat}$ from $\mathrm{Pth}_{\boldsymbol{\mathcal{A}}^{(2)}}$ to $\mathrm{Pth}_{\boldsymbol{\mathcal{B}}^{(2)},\varphi}$ is not necessarily a $\Sigma^{\boldsymbol{\mathcal{A}}^{(2)}}$-homomorphism from $\mathbf{Pth}_{\boldsymbol{\mathcal{A}}^{(2)}}$ to $\mathbf{Pth}_{\boldsymbol{\mathcal{B}}^{(2)}}^{\mathbf{f}^{(2)}}$. 

For the set of sorts $S=1$, the  signature $\Sigma=\{\sigma\}$, where $\sigma$ is a binary operation, and for $X=\{x,y,z\}$, we consider the free $\Sigma$-algebra on $X$, that is $\mathbf{T}_{\Sigma}(X)$. In the example, we will omit any reference to the unique sort in $S$. Consider the set of rewrite rules $\mathcal{A}=\{\mathfrak{p},\mathfrak{q},\mathfrak{r}\}$, where
\begin{align*}
\mathfrak{p}&=(x,y),&
\mathfrak{q}&=(y,z),&
\mathfrak{r}&=(x,z).
\end{align*}
and the set of second-order rewrite rules $\mathcal{A}^{(2)}=\{\mathfrak{p}^{(2)}\}$, where 
$$\mathfrak{p}^{(2)}=([\mathfrak{q}\circ^{0}\mathfrak{p}],[\mathfrak{r}]).$$

Let us note that $\mathfrak{p}^{(2)}$ is a valid second-order rewrite rule because
\allowdisplaybreaks
\begin{align*}
\mathrm{sc}^{(0,1)}(\mathrm{ip}^{(1,X)@}(\mathfrak{q}\circ^{0}\mathfrak{p}))
&=
\mathrm{sc}^{(0,1)}(\mathrm{ip}^{(1,X)@}(\mathfrak{r}))
=
x;
\\
\mathrm{tg}^{(0,1)}(\mathrm{ip}^{(1,X)@}(\mathfrak{q}\circ^{0}\mathfrak{p}))
&=
\mathrm{tg}^{(0,1)}(\mathrm{ip}^{(1,X)@}(\mathfrak{r}))
=z.
\end{align*}

Consider the second-order paths
$$
\xymatrix@C=150pt@R=15pt{
\mathfrak{P}^{(2)}: [
\sigma(\mathfrak{q} \circ^{0} \mathfrak{p}, \mathfrak{q} \circ^{0} \mathfrak{p})
]
\ar@2{->}[r]^-{(\mathfrak{p}^{(2)}, \sigma(\mathfrak{q}\circ^{0}\mathfrak{p}, \underline{\quad}))}
&
[
\sigma(\mathfrak{q}\circ^{0}\mathfrak{p}, \mathfrak{r})
];
\\
\quad
\mathfrak{Q}^{(2)}: [
\sigma(\mathfrak{q}\circ^{0}\mathfrak{p}, \mathfrak{r})
]
\quad
\ar@2{->}[r]^-{(\mathfrak{p}^{(2)}, \sigma(\underline{\quad}, \mathfrak{r}))}
&
\quad
[
\sigma(\mathfrak{r}, \mathfrak{r})
].
\quad
}
$$

One can easily verify that  $\mathfrak{P}^{(2)}$ is a second-order path from $[\sigma(\mathfrak{q}\circ^{0}\mathfrak{p}, \mathfrak{q}\circ^{0}\mathfrak{p})]$ to $[\sigma(\mathfrak{q}\circ^{0}\mathfrak{p},\mathfrak{r})]$ of length $1$ and that $\mathfrak{Q}^{(2)}$ is a second-order path from $[\sigma(\mathfrak{q}\circ^{0}\mathfrak{p}, \mathfrak{r})]$ to $[\sigma(\mathfrak{r},\mathfrak{r})]$ of length $1$. Let us note that none of the one-step subpaths of $\mathfrak{P}^{(2)}$ or $\mathfrak{Q}^{(2)}$ is a  second-order echelon. Thus $\mathfrak{P}^{(2)}$ and $\mathfrak{Q}^{(2)}$ are echelonless second-order paths. Moreover, since both of them are paths of length $1$, it is straightforward to verify that they are both coherent. Finally, both of them are head-constant paths associated with the operation symbol $\sigma$.

Consider the second-order morphism from $\boldsymbol{\mathcal{A}}^{(2)}$ to $\boldsymbol{\mathcal{A}}^{(2)}$ defined as follows:
\begin{enumerate}
\item
Its underlying first-order morphism is $\mathrm{id}^{\boldsymbol{\mathcal{A}}^{(1)}}$, the identity first-order identity mapping at $\boldsymbol{\mathcal{A}}^{(1)}$, introduced in Definition~\ref{DIdRws1};
\item
the second-order echelon mapping $\mathrm{ech}_{\boldsymbol{\mathcal{A}}^{(2)}}^{(2,\mathcal{A}^{(2)})}$ from $\mathcal{A}^{(2)}$ to $\mathrm{Pth}_{\boldsymbol{\mathcal{A}}^{(2)}}$, introduced in Definition~\ref{DDEch}.
\end{enumerate}

We will show that
$$
\mathrm{ech}_{\boldsymbol{\mathcal{A}}^{(2)}}^{(2,\mathcal{A}^{(2)})\flat} \left(
\mathfrak{Q}^{(2)} \circ^{1\mathbf{Pth}_{\boldsymbol{\mathcal{A}}^{(2)}}} \mathfrak{P}^{(2)}
\right)
\neq
\mathrm{ech}_{\boldsymbol{\mathcal{A}}^{(2)}}^{(2,\mathcal{A}^{(2)})\flat} \left(
\mathfrak{Q}^{(2)}
\right)
\circ^{1\mathbf{Pth}_{\boldsymbol{\mathcal{A}}^{(2)}}}
\mathrm{ech}_{\boldsymbol{\mathcal{A}}^{(2)}}^{(2,\mathcal{A}^{(2)})\flat} \left(
\mathfrak{P}^{(2)}
\right).
$$

Consider the right hand side of the equation. Let us recall that $\mathfrak{P}^{(2)}$ is a head-constant echelonless second-order path associated with the operation symbol $\sigma$. Therefore, we are in a position were we can apply the second-order path extraction algorithm from Lemma~\ref{LDPthExtract} which retrieves the following two second-order paths
$$
\xymatrix@C=150pt@R=15pt{
\mathfrak{R}^{(2)}_{0}: [
\mathfrak{q} \circ^{0} \mathfrak{p}
];
\\
\mathfrak{R}^{(2)}_{1}: [
\mathfrak{q}\circ^{0}\mathfrak{p}
]
\ar@2{->}[r]^-{(\mathfrak{p}^{(2)}, \underline{\quad})}
&
[
\mathfrak{r}
].
}
$$
That is, $\mathfrak{R}^{(2)}_{0}$ is the $(2,[1])$-identity path at $[\mathfrak{q} \circ^{0} \mathfrak{p}]$ and $\mathfrak{R}^{(2)}_{1}$ is the second-order echelon associated with the second-order rule $\mathfrak{p}^{(2)}$. Thus, following the definition of the second-order path extension mapping, introduced in Proposition~\ref{PDPthExt},
\begin{align*}
\mathrm{ech}_{\boldsymbol{\mathcal{A}}^{(2)}}^{(2,\mathcal{A}^{(2)})\flat}\left(
\mathfrak{R}^{(2)}_{0}
\right)
&=
\mathrm{ech}_{\boldsymbol{\mathcal{A}}^{(2)}}^{(2,\mathcal{A}^{(2)})\flat}\left(
\mathrm{ip}_{\boldsymbol{\mathcal{A}}^{(2)}}^{(2,[1])\sharp}\left(
[\mathfrak{q} \circ^{0} \mathfrak{p}]
\right)
\right)
=
\mathrm{ip}_{\boldsymbol{\mathcal{A}}^{(2)}}^{(2,[1])\sharp}\left(
[\mathfrak{q} \circ^{0} \mathfrak{p}]
\right)
=
\mathfrak{R}^{(2)}_{0};
\\
\mathrm{ech}_{\boldsymbol{\mathcal{A}}^{(2)}}^{(2,\mathcal{A}^{(2)})\flat}\left(
\mathfrak{R}^{(2)}_{1}
\right)
&=
\mathrm{ech}_{\boldsymbol{\mathcal{A}}^{(2)}}^{(2,\mathcal{A}^{(2)})\flat}\left(
\mathrm{ech}_{\boldsymbol{\mathcal{A}}^{(2)}}^{(2,\mathcal{A}^{(2)})}\left(
\mathfrak{p}^{(2)}
\right)
\right)
=
\mathrm{ech}_{\boldsymbol{\mathcal{A}}^{(2)}}^{(2,\mathcal{A}^{(2)})}\left(
\mathfrak{p}^{(2)}
\right)
=
\mathfrak{R}^{(2)}_{1}.
\end{align*}
Moreover, since $\mathfrak{P}^{(2)}$ is a coherent head-constant echelonless second-order path associated with the operation symbol $\sigma$, following the definition of the second-order path extension mapping, introduced in Proposition~\ref{PDPthExt},
\begin{align*}
\mathrm{ech}_{\boldsymbol{\mathcal{A}}^{(2)}}^{(2,\mathcal{A}^{(2)})\flat}\left(
\mathfrak{P}^{(2)}
\right)
&=
\sigma^{\mathbf{Pth}_{\boldsymbol{\mathcal{A}}^{(2)}}}\left(
\left(
\mathrm{ech}_{\boldsymbol{\mathcal{A}}^{(2)}}^{(2,\mathcal{A}^{(2)})\flat}\left(
\mathfrak{R}^{(2)}_{i}
\right)
\right)_{i \in 2}
\right)
\\
&=
\sigma^{\mathbf{Pth}_{\boldsymbol{\mathcal{A}}^{(2)}}}\left(
\left(
\mathfrak{R}^{(2)}_{i}
\right)_{i \in 2}
\right)
\\
&=
\mathfrak{P}^{(2)}.
\end{align*}

Following a similar argument, it follows that $\mathrm{ech}_{\boldsymbol{\mathcal{A}}^{(2)}}^{(2,\mathcal{A}^{(2)})\flat}(\mathfrak{Q}^{(2)})$ is equal to $\mathfrak{Q}^{(2)}$. Thus, the $1$-composition
$$
\mathrm{ech}_{\boldsymbol{\mathcal{A}}^{(2)}}^{(2,\mathcal{A}^{(2)})\flat} \left(
\mathfrak{Q}^{(2)}
\right)
\circ^{1\mathbf{Pth}_{\boldsymbol{\mathcal{A}}^{(2)}}}
\mathrm{ech}_{\boldsymbol{\mathcal{A}}^{(2)}}^{(2,\mathcal{A}^{(2)})\flat} \left(
\mathfrak{P}^{(2)}
\right)
$$
simplifies to $\mathfrak{Q}^{(2)} \circ^{1\mathbf{Pth}_{\boldsymbol{\mathcal{A}}^{(2)}}} \mathfrak{P}^{(2)}$ which is the second-order path
$$
\xymatrix@C=150pt@R=15pt{
[
\sigma(\mathfrak{q} \circ^{0} \mathfrak{p}, \mathfrak{q} \circ^{0} \mathfrak{p})
]
\ar@2{->}[r]^-{(\mathfrak{p}^{(2)}, \sigma(\mathfrak{q}\circ^{0}\mathfrak{p}, \underline{\quad}))}
&
[
\sigma(\mathfrak{q}\circ^{0}\mathfrak{p}, \mathfrak{r})
];
\\
{\quad\quad\quad\quad\quad\quad\quad}
\ar@2{->}[r]^-{(\mathfrak{p}^{(2)}, \sigma(\underline{\quad}, \mathfrak{r}))}
&
\quad
[
\sigma(\mathfrak{r}, \mathfrak{r})
].
\quad
}
$$

Now consider the left hand side of the equation. One can easily verify that the $1$-composition $\mathfrak{Q}^{(2)} \circ^{1\mathbf{Pth}_{\boldsymbol{\mathcal{A}}^{(2)}}} \mathfrak{P}^{(2)}$ is an echelonless, coherent, head-constant second-order path associated with the operation symbol $\sigma$. Therefore, we are in a position were we can apply the second-order path extraction algorithm from Lemma~\ref{LDPthExtract} which retrieves the following two second-order paths
$$
\xymatrix@C=150pt@R=15pt{
\mathfrak{R}^{(2)\prime}_{0}: [
\mathfrak{q} \circ^{0} \mathfrak{p}
]
\ar@2{->}[r]^-{(\mathfrak{p}^{(2)}, \underline{\quad})}
&
[
\mathfrak{r}
];
\\
\mathfrak{R}^{(2)\prime}_{1}: [
\mathfrak{q}\circ^{0}\mathfrak{p}
]
\ar@2{->}[r]^-{(\mathfrak{p}^{(2)}, \underline{\quad})}
&
[
\mathfrak{r}
].
}
$$
That is, they both are the second-order echelon associated with the second-order rule $\mathfrak{p}^{(2)}$. Thus, following the definition of the second-order path extension mapping, introduced in Proposition~\ref{PDPthExt},
\begin{align*}
\mathrm{ech}_{\boldsymbol{\mathcal{A}}^{(2)}}^{(2,\mathcal{A}^{(2)})\flat}\left(
\mathfrak{R}^{(2)\prime}_{0}
\right)
&=
\mathrm{ech}_{\boldsymbol{\mathcal{A}}^{(2)}}^{(2,\mathcal{A}^{(2)})\flat}\left(
\mathrm{ech}_{\boldsymbol{\mathcal{A}}^{(2)}}^{(2,\mathcal{A}^{(2)})}\left(
\mathfrak{p}^{(2)}
\right)
\right)
=
\mathrm{ech}_{\boldsymbol{\mathcal{A}}^{(2)}}^{(2,\mathcal{A}^{(2)})}\left(
\mathfrak{p}^{(2)}
\right)
=
\mathfrak{R}^{(2)\prime}_{0};
\\
\mathrm{ech}_{\boldsymbol{\mathcal{A}}^{(2)}}^{(2,\mathcal{A}^{(2)})\flat}\left(
\mathfrak{R}^{(2)\prime}_{1}
\right)
&=
\mathrm{ech}_{\boldsymbol{\mathcal{A}}^{(2)}}^{(2,\mathcal{A}^{(2)})\flat}\left(
\mathrm{ech}_{\boldsymbol{\mathcal{A}}^{(2)}}^{(2,\mathcal{A}^{(2)})}\left(
\mathfrak{p}^{(2)}
\right)
\right)
=
\mathrm{ech}_{\boldsymbol{\mathcal{A}}^{(2)}}^{(2,\mathcal{A}^{(2)})}\left(
\mathfrak{p}^{(2)}
\right)
=
\mathfrak{R}^{(2)\prime}_{1}.
\end{align*}
Moreover, since $\mathfrak{Q}^{(2)} \circ^{1\mathbf{Pth}_{\boldsymbol{\mathcal{A}}^{(2)}}} \mathfrak{P}^{(2)}$ is a coherent head-constant echelonless second-order path associated with the operation symbol $\sigma$, following the definition of the second-order path extension mapping, introduced in Proposition~\ref{PDPthExt},
\begin{align*}
\mathrm{ech}_{\boldsymbol{\mathcal{A}}^{(2)}}^{(2,\mathcal{A}^{(2)})\flat}\left(
\mathfrak{Q}^{(2)} \circ^{1\mathbf{Pth}_{\boldsymbol{\mathcal{A}}^{(2)}}} \mathfrak{P}^{(2)}
\right)
&=
\sigma^{\mathbf{Pth}_{\boldsymbol{\mathcal{A}}^{(2)}}}\left(
\left(
\mathrm{ech}_{\boldsymbol{\mathcal{A}}^{(2)}}^{(2,\mathcal{A}^{(2)})\flat}\left(
\mathfrak{R}^{(2)\prime}_{i}
\right)
\right)_{i \in 2}
\right)
\\
&=
\sigma^{\mathbf{Pth}_{\boldsymbol{\mathcal{A}}^{(2)}}}\left(
\left(
\mathfrak{R}^{(2)\prime}_{i}
\right)_{i \in 2}
\right).
\end{align*}
That is, $\sigma^{\mathbf{Pth}_{\boldsymbol{\mathcal{A}}^{(2)}}}(\mathfrak{R}^{(2)\prime}_{0}, \mathfrak{R}^{(2)\prime}_{1})$ is the second-order path 
$$
\xymatrix@C=150pt@R=15pt{
[
\sigma(\mathfrak{q} \circ^{0} \mathfrak{p}, \mathfrak{q} \circ^{0} \mathfrak{p})
]
\ar@2{->}[r]^-{(\mathfrak{p}^{(2)}, \sigma(\underline{\quad}, \mathfrak{q}\circ^{0}\mathfrak{p}))}
&
[
\sigma(\mathfrak{r}, \mathfrak{q}\circ^{0}\mathfrak{p})
];
\\
{\quad\quad\quad\quad\quad\quad\quad}
\ar@2{->}[r]^-{(\mathfrak{p}^{(2)}, \sigma(\mathfrak{r}, \underline{\quad}))}
&
\quad
[
\sigma(\mathfrak{r}, \mathfrak{r})
]
\quad
}
$$
which is not equal to the second-order path $\mathfrak{Q}^{(2)} \circ^{1\mathbf{Pth}_{\boldsymbol{\mathcal{A}}^{(2)}}} \mathfrak{P}^{(2)}$.

All in all, the second-order path extension mapping does not respect the $1$-composition of second-order paths and, therefore, is not necessarily a $\Sigma^{\boldsymbol{\mathcal{A}}^{(2)}}$-homomorphism.
\end{remark}

Now, given a second-order morphism $\mathbf{f}^{(2)}$, $\llbracket\mathrm{Pth}_{\boldsymbol{\mathcal{B}}^{(2)}}\rrbracket_{\varphi}$ is equipped in a natural way with a structure of partial $\Sigma^{\boldsymbol{\mathcal{A}}^{(2)}}$-algebra that we denote by $\llbracket\mathbf{Pth}_{\boldsymbol{\mathcal{B}}^{(2)}}^{\mathbf{f}^{(2)}}\rrbracket$. To end this section we show that $\mathrm{pr}^{\llbracket\cdot\rrbracket}_{\boldsymbol{\mathcal{B}}^{(2)}, \varphi}$ is a surjective $\Sigma^{\boldsymbol{\mathcal{A}}^{(2)}}$-homomorphism from $\mathbf{Pth}_{\boldsymbol{\mathcal{B}}^{(2)}}^{\mathbf{f}^{(2)}}$ to $\llbracket\mathbf{Pth}_{\boldsymbol{\mathcal{B}}^{(2)}}^{\mathbf{f}^{(2)}}\rrbracket$.

\begin{proposition}
\label{PDQPthBDCatAlg}
Let $\mathbf{f}^{(2)}=(\varphi, c, (f^{(i)})_{i\in 3})$ be a second-order morphism from $\boldsymbol{\mathcal{A}}^{(2)}$ to $\boldsymbol{\mathcal{B}}^{(2)}$. Then the $S$-sorted set $\llbracket\mathrm{Pth}_{\boldsymbol{\mathcal{B}}^{(2)}}\rrbracket_{\varphi}$ is equipped, in a natural way, with a structure of partial $\Sigma^{\boldsymbol{\mathcal{A}}^{(2)}}$-algebra.
\end{proposition}

\begin{proof}
Let us denote by $\llbracket\mathbf{Pth}_{\boldsymbol{\mathcal{B}}^{(2)}}^{\mathbf{f}^{(2)}}\rrbracket$ the $\Sigma^{\boldsymbol{\mathcal{A}}^{(2)}}$-algebra defined as follows:

\textsf{(1)}
The underlying $S$-sorted set of $\llbracket\mathbf{Pth}_{\boldsymbol{\mathcal{B}}^{(2)}}^{\mathbf{f}^{(2)}}\rrbracket$ is $\llbracket\mathrm{Pth}_{\boldsymbol{\mathcal{B}}^{(2)}}\rrbracket_{\varphi}=(\llbracket\mathrm{Pth}_{\boldsymbol{\mathcal{B}}^{(2)}}\rrbracket_{\varphi(s)})_{s\in S}$.

\textsf{(2)}
For every $(\mathbf{s}, s)\in S^{\ast}\times S$ and every operation symbol $\sigma\in\Sigma_{\mathbf{s},s}$, the operation $\sigma^{\llbracket\mathbf{Pth}_{\boldsymbol{\mathcal{B}}^{(2)}}^{\mathbf{f}^{(2)}}\rrbracket}$ is given by the interpretation of $\sigma$ in the $\Sigma$-algebra $\llbracket\mathbf{Pth}_{\boldsymbol{\mathcal{B}}^{(2)}}^{\mathbf{f}^{(2)}(0,2)}\rrbracket$ introduced in Remark~\ref{RDDSigmaAlg}.

\textsf{(3)}
For every $s\in S$ and every $\mathfrak{p}\in\mathcal{A}_{s}^{(1)}$, the constant $\mathfrak{p}^{\llbracket\mathbf{Pth}_{\boldsymbol{\mathcal{B}}^{(2)}}^{\mathbf{f}^{(2)}}\rrbracket}$ is given by the interpretation of $\mathfrak{p}$ in $\llbracket\mathbf{Pth}_{\boldsymbol{\mathcal{A}}^{(2)}}^{\mathbf{f}^{(2)}(1,2)}\rrbracket$ that, we recall, was stated in Proposition~\ref{PDQPthBCatAlg}.

\textsf{(4)}
For every $s\in S$, the interpretations of the operations $\mathrm{sc}_{s}^{0}$ and $\mathrm{tg}_{s}^{0}$ are given by the interpretation of $\mathrm{sc}_{s}^{0}$ and $\mathrm{tg}_{s}^{0}$ in $\llbracket\mathbf{Pth}_{\boldsymbol{\mathcal{A}}^{(2)}}^{\mathbf{f}^{(2)}(1,2)}\rrbracket$ that , we recall, was stated in Proposition~\ref{PDQPthBCatAlg}.

\textsf{(5)}
Similarly,  for every $s \in S$, the interpretations of the operation $\circ_{s}^{0}$ is given by the interpretation of $\circ_{s}^{0}$ in $\llbracket\mathbf{Pth}_{\boldsymbol{\mathcal{A}}^{(2)}}^{\mathbf{f}^{(2)}(1,2)}\rrbracket$ that , we recall, was stated in Proposition~\ref{PDQPthBCatAlg}.

\textsf{(6)}
For every $s\in S$ and every $\mathfrak{p}^{(2)}\in\mathcal{A}_{s}^{(2)}$, the constant $\mathfrak{p}^{(2)\llbracket\mathbf{Pth}_{\boldsymbol{\mathcal{B}}^{(2)}}^{\mathbf{f}^{(2)}}\rrbracket}$ is given by the class of the image of $\mathfrak{p}^{(2)}$ under the mapping $f^{(2)}$, i.e., we set
$$
\mathfrak{p}^{(2)\llbracket\mathbf{Pth}_{\boldsymbol{\mathcal{B}}^{(2)}}^{\mathbf{f}^{(2)}}\rrbracket}
=
\left\llbracket
f^{(2)}_{s}\left(
\mathfrak{p}^{(2)}
\right)
\right\rrbracket_{\varphi(s)}
\in
\llbracket
\mathrm{Pth}_{\boldsymbol{\mathcal{B}}^{(2)}}
\rrbracket_{\varphi(s)}
$$.

\textsf{(7)}
For every $s\in S$, the interpretations of the operations $\mathrm{sc}_{s}^{1}$ and $\mathrm{tg}_{s}^{1}$ are given by
\allowdisplaybreaks
\begin{align*}
\mathrm{sc}_{s}^{1\llbracket\mathbf{Pth}_{\boldsymbol{\mathcal{B}}^{(2)}}^{\mathbf{f}^{(2)}}\rrbracket}
&=
\mathrm{sc}_{\varphi(s)}^{1\llbracket\mathbf{Pth}_{\boldsymbol{\mathcal{B}}^{(2)}}\rrbracket}
&&\mbox{and}&
\mathrm{tg}_{s}^{1\llbracket\mathbf{Pth}_{\boldsymbol{\mathcal{B}}^{(2)}}^{\mathbf{f}^{(2)}}\rrbracket}
&=
\mathrm{tg}_{\varphi(s)}^{1\llbracket\mathbf{Pth}_{\boldsymbol{\mathcal{B}}^{(2)}}\rrbracket}.
\end{align*}
That is, their interpretations are given by the interpretations of $\mathrm{sc}_{\varphi(s)}^{1}$ and $\mathrm{tg}_{\varphi(s)}^{1}$ in $\llbracket\mathbf{Pth}_{\boldsymbol{\mathcal{B}}^{(2)}}\rrbracket$ that, we recall, was introduced in Proposition~\ref{PDVDCatAlg}.

\textsf{(8)}
Similarly, for every $s \in S$,  the interpretation of the partial binary operation $\circ_{s}^{1}$ is given by
$$
\circ_{s}^{1\llbracket\mathbf{Pth}_{\boldsymbol{\mathcal{B}}^{(2)}}^{\mathbf{f}^{(2)}}\rrbracket}
=
\circ_{\varphi(s)}^{1\llbracket\mathbf{Pth}_{\boldsymbol{\mathcal{B}}^{(2)}}\rrbracket}.
$$
That is, its interpretation is given by the interpretation of $\circ_{\varphi(s)}^{1}$ in $\llbracket\mathbf{Pth}_{\boldsymbol{\mathcal{B}}^{(2)}}\rrbracket$ that, we recall, was introduced in Proposition~\ref{PDVDCatAlg}.

This completes the definition of the partial $\Sigma^{\boldsymbol{\mathcal{A}}^{(2)}}$-algebra $\llbracket\mathbf{Pth}_{\boldsymbol{\mathcal{B}}^{(2)}}^{\mathbf{f}^{(2)}}\rrbracket$.
\end{proof}

\begin{proposition}
\label{PDKerCHDCatHom}
Let $\mathbf{f}^{(2)}=(\varphi, c, (f^{(i)})_{i\in 3})$ be a second-order morphism from $\boldsymbol{\mathcal{A}}^{(2)}$ to $\boldsymbol{\mathcal{B}}^{(2)}$. Then the mapping 
$$
\mathrm{pr}_{\boldsymbol{\mathcal{B}}^{(2)}, \varphi}^{\llbracket \cdot \rrbracket}
\colon 
\mathrm{Pth}_{\boldsymbol{\mathcal{B}}^{(2)}, \varphi}
\mor
\llbracket\mathrm{Pth}_{\boldsymbol{\mathcal{B}}^{(2)}}\rrbracket_{\varphi}
$$
is a surjective $\Sigma^{\boldsymbol{\mathcal{A}}^{(2)}}$-homomorphism from $\mathbf{Pth}_{\boldsymbol{\mathcal{B}}^{(2)}}^{\mathbf{f}^{(2)}}$ to $\llbracket\mathbf{Pth}_{\boldsymbol{\mathcal{B}}^{(2)}}^{\mathbf{f}^{(2)}}\rrbracket$.
\end{proposition}

\begin{proof}
We prove that $\mathrm{pr}_{\boldsymbol{\mathcal{B}}^{(2)}, \varphi}^{\llbracket \cdot \rrbracket}$ is compatible with every operation symbol in $\Sigma^{\boldsymbol{\mathcal{A}}^{(2)}}$.

{\sffamily The mapping $\mathrm{pr}_{\boldsymbol{\mathcal{B}}^{(2)}, \varphi}^{\llbracket \cdot \rrbracket}$ is a $\Sigma^{\boldsymbol{\mathcal{A}}^{(1)}}$-homomorphism.}

The fact that the mapping $\mathrm{pr}_{\boldsymbol{\mathcal{B}}^{(2)}, \varphi}^{\llbracket \cdot \rrbracket}$ is a $\Sigma^{\boldsymbol{\mathcal{A}}^{(1)}}$-homomorphism follows from Proposition~\ref{PDPrBCatHom}.

{\sffamily The mapping $\mathrm{pr}_{\boldsymbol{\mathcal{B}}^{(2)}, \varphi}^{\llbracket \cdot \rrbracket}$ is compatible with the second-order rewrite rules.}

Let $s$ be a sort in $S$ and $\mathfrak{p}^{(2)}$ a rewrite rule in $\mathcal{A}_{s}^{(2)}$. Thus,
$$
\left\llbracket
\mathfrak{p}^{(2)\mathbf{Pth}_{\boldsymbol{\mathcal{B}}^{(2)}}^{\mathbf{f}^{(2)}}}
\right\rrbracket_{\varphi(s)}
=
\left\llbracket
f_{s}^{(2)\flat}\left(
\mathfrak{p}^{(2)}
\right)
\right\rrbracket_{\varphi(s)}
=
\mathfrak{p}^{(2)\llbracket\mathbf{Pth}_{\boldsymbol{\mathcal{B}}^{(2)}}^{\mathbf{f}^{(2)}}\rrbracket}.
$$

Hence, $\mathrm{pr}_{\boldsymbol{\mathcal{B}}^{(2)}, \varphi}^{\llbracket \cdot \rrbracket}$ is compatible with the second-order rewrite rules.

{\sffamily The mapping $\mathrm{pr}_{\boldsymbol{\mathcal{B}}^{(2)}, \varphi}^{\llbracket \cdot \rrbracket}$ is compatible with the $1$-source.}

Let $s$ be a sort in $S$ and let us consider the $1$-source operation symbol $\mathrm{sc}_{s}^{1}$ in $\Sigma^{\boldsymbol{\mathcal{A}}^{(2)}}_{s,s}$. Let $\mathfrak{P}^{(2)}$ be a second-order path in $\mathrm{Pth}_{\boldsymbol{\mathcal{B}}^{(2)}, \varphi(s)}$.

The following chain of equalities holds
\allowdisplaybreaks
\begin{align*}
\left\llbracket
\mathrm{sc}_{s}^{1\mathbf{Pth}_{\boldsymbol{\mathcal{B}}^{(2)}}^{\mathbf{f}^{(2)}}}\left(
\mathfrak{P}^{(2)}
\right)
\right\rrbracket_{\varphi(s)}
&=
\left\llbracket
\mathrm{sc}_{\varphi(s)}^{1\mathbf{Pth}_{\boldsymbol{\mathcal{B}}^{(2)}}}\left(
\mathfrak{P}^{(2)}
\right)
\right\rrbracket_{\varphi(s)}
\tag{1}
\\
&=
\mathrm{sc}_{\varphi(s)}^{1\llbracket\mathbf{Pth}_{\boldsymbol{\mathcal{B}}^{(2)}}\rrbracket}\left(
\left\llbracket
\mathfrak{P}^{(2)}
\right\rrbracket_{\varphi(s)}
\right)
\tag{2}
\\
&=
\mathrm{sc}_{s}^{1\llbracket\mathbf{Pth}_{\boldsymbol{\mathcal{B}}^{(2)}}^{\mathbf{f}^{(2)}}\rrbracket}\left(
\left\llbracket
\mathfrak{P}^{(2)}
\right\rrbracket_{\varphi(s)}
\right).
\tag{3}
\end{align*}
The first equality unravels the interpretation of the operation symbol $\mathrm{sc}_{s}^{1}$ in the partial $\Sigma^{\boldsymbol{\mathcal{A}}^{(2)}}$-algebra $\mathbf{Pth}_{\boldsymbol{\mathcal{B}}^{(2)}}^{\mathbf{f}^{(2)}}$, introduced in Proposition~\ref{PDPthBDCatAlg};
the second equality follows from the fact that, according to Proposition~\ref{PDVDCatAlg}, $\mathrm{pr}^{\llbracket\cdot\rrbracket}_{\boldsymbol{\mathcal{B}}^{(2)}}$ is a $\Lambda^{\boldsymbol{\mathcal{B}}^{(2)}}$-homomorphism from $\mathbf{Pth}_{\boldsymbol{\mathcal{B}}^{(2)}}$ to $\llbracket\mathbf{Pth}_{\boldsymbol{\mathcal{B}}^{(2)}}\rrbracket$;
finally, the last equality recovers the interpretation of the operation symbol $\mathrm{sc}_{s}^{1}$ in the partial $\Sigma^{\boldsymbol{\mathcal{A}}^{(2)}}$-algebra $\llbracket\mathbf{Pth}_{\boldsymbol{\mathcal{B}}^{(2)}}^{\mathbf{f}^{(2)}}\rrbracket$, introduced in Proposition~\ref{PDQPthBDCatAlg}.

Hence, $\mathrm{pr}_{\boldsymbol{\mathcal{B}}^{(2)}, \varphi}^{\llbracket \cdot \rrbracket}$ is compatible with the $1$-source operation.

{\sffamily The mapping $\mathrm{pr}_{\boldsymbol{\mathcal{B}}^{(2)}, \varphi}^{\llbracket \cdot \rrbracket}$ is compatible with the $1$-target.}

Let $s$ be a sort in $S$ and let us consider the $1$-source operation symbol $\mathrm{tg}_{s}^{1}$ in $\Sigma^{\boldsymbol{\mathcal{A}}^{(2)}}_{s,s}$. Let $\mathfrak{P}^{(2)}$ be a second-order path in $\mathrm{Pth}_{\boldsymbol{\mathcal{B}}^{(2)}, \varphi(s)}$, then the following equality holds
$$
\left\llbracket
\mathrm{tg}_{s}^{1\mathbf{Pth}_{\boldsymbol{\mathcal{B}}^{(2)}}^{\mathbf{f}^{(2)}}}\left(
\mathfrak{P}^{(2)}
\right)
\right\rrbracket_{\varphi(s)}
=
\mathrm{tg}_{s}^{1\llbracket\mathbf{Pth}_{\boldsymbol{\mathcal{B}}^{(2)}}^{\mathbf{f}^{(2)}}\rrbracket}\left(
\left\llbracket
\mathfrak{P}^{(2)}
\right\rrbracket_{\varphi(s)}
\right).
$$

The proof of this case is identical to that of the $1$-source.

Hence, $\mathrm{pr}_{\boldsymbol{\mathcal{B}}^{(2)}, \varphi}^{\llbracket \cdot \rrbracket}$ is compatible with the $1$-target operation.

{\sffamily The mapping $\mathrm{pr}_{\boldsymbol{\mathcal{B}}^{(2)}, \varphi}^{\llbracket \cdot \rrbracket}$ is compatible with the $1$-composition.}

Let $s$ be a sort in $S$ and let us consider the $1$-composition operation symbol $\circ_{s}^{1}$ in $\Sigma_{ss,s}^{\boldsymbol{\mathcal{A}}^{(2)}}$. Let $\mathfrak{P}^{(2)}$ and $\mathfrak{Q}^{(2)}$ be two second-order paths in $\mathrm{Pth}_{\boldsymbol{\mathcal{B}}^{(2)}, \varphi(s)}$ such that
$$
\mathrm{sc}_{\boldsymbol{\mathcal{B}}^{(2)}, \varphi(s)}^{([1],2)}\left(\mathfrak{Q}^{(2)}\right)
=
\mathrm{tg}_{\boldsymbol{\mathcal{B}}^{(2)}, \varphi(s)}^{([1],2)}\left(\mathfrak{P}^{(2)}\right).
$$
Then the following equality holds
$$
\left\llbracket
\mathfrak{Q}^{(2)}
\circ_{s}^{1\mathbf{Pth}_{\boldsymbol{\mathcal{B}}^{(2)}}^{\mathbf{f}^{(2)}}}
\mathfrak{P}^{(2)}
\right\rrbracket_{\varphi(s)}
=
\left\llbracket
\mathfrak{Q}^{(2)}
\right\rrbracket_{\varphi(s)}
\circ_{s}^{1\llbracket\mathbf{Pth}_{\boldsymbol{\mathcal{B}}^{(2)}}^{\mathbf{f}^{(2)}}\rrbracket}
\left\llbracket
\mathfrak{P}^{(2)}
\right\rrbracket_{\varphi(s)}.
$$

The proof of this case is identical to that of the $1$-source.

Hence, $\mathrm{pr}_{\boldsymbol{\mathcal{B}}^{(2)}, \varphi}^{\llbracket \cdot \rrbracket}$ is compatible with the $1$-composition operation.

This completes the proof.
\end{proof}

\section{A structure of partial $\Sigma^{\boldsymbol{\mathcal{A}}^{(2)}}$-algebra on $\mathrm{PT}_{\boldsymbol{\mathcal{B}}^{(2)}, \varphi}$}

In this section we show that, given a second-order morphism $\mathbf{f}^{(2)}$, $\mathrm{PT}_{\boldsymbol{\mathcal{B}}^{(2)}, \varphi}$ is equipped in a natural way with a structure of partial $\Sigma^{\boldsymbol{\mathcal{A}}^{(2)}}$-algebra that we denote by $\mathbf{PT}_{\boldsymbol{\mathcal{B}}^{(2)}}^{\mathbf{f}^{(2)}}$.

\begin{proposition}
\label{PDPTBDCatAlg}
Let $\mathbf{f}^{(2)}=(\varphi, c, (f^{(i)})_{i\in 3})$ be a second-order morphism from $\boldsymbol{\mathcal{A}}^{(2)}$ to $\boldsymbol{\mathcal{B}}^{(2)}$. Then the $S$-sorted set $\mathrm{PT}_{\boldsymbol{\mathcal{B}}^{(2)}, \varphi}$ is equipped, in a natural way, with a structure of partial $\Sigma^{\boldsymbol{\mathcal{A}}^{(2)}}$-algebra.
\end{proposition}

\begin{proof}
Let us denote by $\mathbf{PT}_{\boldsymbol{\mathcal{B}}^{(2)}}^{\mathbf{f}^{(2)}}$ the $\Sigma^{\boldsymbol{\mathcal{A}}^{(2)}}$-algebra defined as follows:

\textsf{(1)}
The underlying $S$-sorted set of $\mathbf{PT}_{\boldsymbol{\mathcal{B}}^{(2)}}^{\mathbf{f}^{(2)}}$ is $\mathrm{PT}_{\boldsymbol{\mathcal{B}}^{(2)}, \varphi}=(\mathrm{PT}_{\boldsymbol{\mathcal{B}}^{(2)}, \varphi(s)})_{s\in S}$.

\textsf{(2)}
For every $(\mathbf{s}, s)\in S^{\ast}\times S$ and every operation symbol $\sigma\in\Sigma_{\mathbf{s},s}$, the operation $\sigma^{\mathbf{PT}_{\boldsymbol{\mathcal{B}}^{(2)}}^{\mathbf{f}^{(2)}}}$ is given by the interpretation of $\sigma$ in the $\Sigma$-algebra $\mathbf{c}_{\mathfrak{d}}^{\ast}(\mathbf{PT}_{\boldsymbol{\mathcal{B}}^{(2)}}^{(0,2)})$. That is, its interpretation is given by the derived operation in $\mathbf{PT}^{(0,2)}_{\boldsymbol{\mathcal{B}}^{(2)}}$, the $\Lambda$-reduct of the partial $\Lambda^{\boldsymbol{\mathcal{B}}^{(2)}}$-algebra $\mathbf{PT}_{\boldsymbol{\mathcal{B}}^{(2)}}$, introduced in Proposition~\ref{PDPTDCatAlg}.

\textsf{(3)}
For every $s\in S$ and every $\mathfrak{p}\in\mathcal{A}_{s}^{(1)}$, the constant $\mathfrak{p}^{\mathbf{PT}_{\boldsymbol{\mathcal{B}}^{(2)}}^{\mathbf{f}^{(2)}}}$ is given by
$$
\mathfrak{p}^{\mathbf{PT}_{\boldsymbol{\mathcal{B}}^{(2)}}^{\mathbf{f}^{(2)}}}
=
\mathrm{CH}^{(2)}_{\boldsymbol{\mathcal{B}}^{(2)}, \varphi(s)}\left(
f^{(2)\flat}_{s}\left(
\mathfrak{p}^{\mathbf{Pth}_{\boldsymbol{\mathcal{A}}^{(2)}}}
\right)
\right).
$$

\textsf{(4)}
For every $s\in S$, the interpretations of the operations $\mathrm{sc}_{s}^{0}$ and $\mathrm{tg}_{s}^{0}$ are given by
\begin{align*}
\mathrm{sc}_{s}^{0\mathbf{PT}_{\boldsymbol{\mathcal{B}}^{(2)}}^{\mathbf{f}^{(2)}}}
&=
\mathrm{sc}_{\varphi(s)}^{0\mathbf{PT}_{\boldsymbol{\mathcal{B}}^{(2)}}}
&&\mbox{and}&
\mathrm{tg}_{s}^{0\mathbf{PT}_{\boldsymbol{\mathcal{B}}^{(2)}}^{\mathbf{f}^{(2)}}}
&=
\mathrm{tg}_{\varphi(s)}^{0\mathbf{PT}_{\boldsymbol{\mathcal{B}}^{(2)}}}.
\end{align*}
That is, their interpretations are given by the interpretations of $\mathrm{sc}_{\varphi(s)}^{0}$ and $\mathrm{tg}_{\varphi(s)}^{0}$ in $\mathbf{PT}_{\boldsymbol{\mathcal{B}}^{(2)}}$ that, we recall, was introduced in Proposition~\ref{PDPTDCatAlg}.

\textsf{(5)}
Similarly,  for every $s \in S$, the interpretations of the operation $\circ_{s}^{0}$ is given by
$$
\circ_{s}^{0\mathbf{PT}_{\boldsymbol{\mathcal{B}}^{(2)}}^{\mathbf{f}^{(2)}}}
=
\circ_{\varphi(s)}^{0\mathbf{PT}_{\boldsymbol{\mathcal{B}}^{(2)}}}.
$$
That is, its interpretation is given by the interpretation of $\circ_{\varphi(s)}^{0}$ in $\mathbf{PT}_{\boldsymbol{\mathcal{B}}^{(2)}}$ that, we recall, was introduced in Proposition~\ref{PDPTDCatAlg}.

\textsf{(6)}
For every $s\in S$ and every $\mathfrak{p}^{(2)}\in\mathcal{A}_{s}^{(2)}$, the constant $\mathfrak{p}^{\mathbf{Pth}_{\boldsymbol{\mathcal{B}}^{(2)}}^{\mathbf{f}^{(2)}}}$ is given by
$$
\mathfrak{p}^{(2)\mathbf{PT}_{\boldsymbol{\mathcal{B}}^{(2)}}^{\mathbf{f}^{(2)}}}
=
\mathrm{CH}^{(2)}_{\boldsymbol{\mathcal{B}}^{(2)}, \varphi(s)}\left(
f^{(2)}_{s}\left(
\mathfrak{p}^{(2)}
\right)
\right)
\in
\mathrm{PT}_{\boldsymbol{\mathcal{B}}^{(2)}, \varphi(s)}.
$$

\textsf{(7)}
For every $s\in S$, the interpretations of the operations $\mathrm{sc}_{s}^{1}$ and $\mathrm{tg}_{s}^{1}$ are given by
\allowdisplaybreaks
\begin{align*}
\mathrm{sc}_{s}^{1\mathbf{PT}_{\boldsymbol{\mathcal{B}}^{(2)}}^{\mathbf{f}^{(2)}}}
&=
\mathrm{sc}_{\varphi(s)}^{1\mathbf{PT}_{\boldsymbol{\mathcal{B}}^{(2)}}}
&&\mbox{and}&
\mathrm{tg}_{s}^{1\mathbf{PT}_{\boldsymbol{\mathcal{B}}^{(2)}}^{\mathbf{f}^{(2)}}}
&=
\mathrm{tg}_{\varphi(s)}^{1\mathbf{PT}_{\boldsymbol{\mathcal{B}}^{(2)}}}.
\end{align*}
That is, their interpretations are given by the interpretations of $\mathrm{sc}_{\varphi(s)}^{1}$ and $\mathrm{tg}_{\varphi(s)}^{1}$ in $\mathbf{PT}_{\boldsymbol{\mathcal{B}}^{(2)}}$ that, we recall, was introduced in Proposition~\ref{PDPTDCatAlg}.

\textsf{(8)}
Similarly, for every $s \in S$,  the interpretation of the partial binary operation $\circ_{s}^{1}$ is given by
$$
\circ_{s}^{1\mathbf{PT}_{\boldsymbol{\mathcal{B}}^{(2)}}^{\mathbf{f}^{(2)}}}
=
\circ_{\varphi(s)}^{1\mathbf{PT}_{\boldsymbol{\mathcal{B}}^{(2)}}}.
$$
That is, its interpretation is given by the interpretation of $\circ_{\varphi(s)}^{1}$ in $\mathbf{PT}_{\boldsymbol{\mathcal{B}}^{(2)}}$ that, we recall, was introduced in Proposition~\ref{PDPTDCatAlg}.

This completes the definition of the partial $\Sigma^{\boldsymbol{\mathcal{A}}^{(2)}}$-algebra $\mathbf{PT}_{\boldsymbol{\mathcal{B}}^{(2)}}^{\mathbf{f}^{(2)}}$.
\end{proof}

Now, given a second-order morphism $\mathbf{f}^{(2)}$, $\llbracket\mathrm{PT}_{\boldsymbol{\mathcal{B}}^{(2)}}\rrbracket_{\varphi}$ is equipped in a natural way with a structure of partial $\Sigma^{\boldsymbol{\mathcal{A}}^{(2)}}$-algebra that we denote by $\llbracket\mathbf{PT}_{\boldsymbol{\mathcal{B}}^{(2)}}^{\mathbf{f}^{(2)}}\rrbracket$. Thus, we show that $\mathrm{pr}^{\Theta^{\llbracket 2 \rrbracket}}_{\boldsymbol{\mathcal{B}}^{(2)}, \varphi}$ is a surjective $\Sigma^{\boldsymbol{\mathcal{A}}^{(2)}}$-homomorphism from $\mathbf{PT}_{\boldsymbol{\mathcal{B}}^{(2)}}^{\mathbf{f}^{(2)}}$ to $\llbracket\mathbf{PT}_{\boldsymbol{\mathcal{B}}^{(2)}}^{\mathbf{f}^{(2)}}\rrbracket$.

\begin{proposition}
\label{PDQPTBDCatAlg}
Let $\mathbf{f}^{(2)}=(\varphi, c, (f^{(i)})_{i\in 3})$ be a second-order morphism from $\boldsymbol{\mathcal{A}}^{(2)}$ to $\boldsymbol{\mathcal{B}}^{(2)}$. Then the $S$-sorted set $\llbracket\mathrm{PT}_{\boldsymbol{\mathcal{B}}^{(2)}}\rrbracket_{\varphi}$ is equipped, in a natural way, with a structure of partial $\Sigma^{\boldsymbol{\mathcal{A}}^{(2)}}$-algebra.
\end{proposition}

\begin{proof}
Let us denote by $\llbracket\mathbf{PT}_{\boldsymbol{\mathcal{B}}^{(2)}}^{\mathbf{f}^{(2)}}\rrbracket$ the $\Sigma^{\boldsymbol{\mathcal{A}}^{(2)}}$-algebra defined as follows:

\textsf{(1)}
The underlying $S$-sorted set of $\llbracket\mathbf{PT}_{\boldsymbol{\mathcal{B}}^{(2)}}^{\mathbf{f}^{(2)}}\rrbracket$ is $\llbracket\mathrm{PT}_{\boldsymbol{\mathcal{B}}^{(2)}}\rrbracket_{\varphi}=(\llbracket\mathrm{PT}_{\boldsymbol{\mathcal{B}}^{(2)}}\rrbracket_{\varphi(s)})_{s\in S}$.

\textsf{(2)}
For every $(\mathbf{s}, s)\in S^{\ast}\times S$ and every operation symbol $\sigma\in\Sigma_{\mathbf{s},s}$, the operation $\sigma^{\mathbf{PT}_{\boldsymbol{\mathcal{B}}^{(2)}}^{\mathbf{f}^{(2)}}}$ is given by the interpretation of $\sigma$ in the $\Sigma$-algebra $\mathbf{c}_{\mathfrak{d}}^{\ast}(\llbracket\mathbf{PT}_{\boldsymbol{\mathcal{B}}^{(2)}}^{(0,2)}\rrbracket)$. That is, its interpretation is given by the derived operation in $\llbracket\mathbf{PT}^{(0,2)}_{\boldsymbol{\mathcal{B}}^{(2)}}\rrbracket$, the $\Lambda$-reduct of the partial $\Lambda^{\boldsymbol{\mathcal{B}}^{(2)}}$-algebra $\llbracket\mathbf{PT}_{\boldsymbol{\mathcal{B}}^{(2)}}\rrbracket$, introduced in Proposition~\ref{PDPTQDCatAlg}.

\textsf{(3)}
For every $s\in S$ and every $\mathfrak{p}\in\mathcal{A}_{s}^{(1)}$, the constant $\mathfrak{p}^{\mathbf{PT}_{\boldsymbol{\mathcal{B}}^{(2)}}^{\mathbf{f}^{(2)}}}$ is given by
$$
\mathfrak{p}^{\llbracket\mathbf{PT}_{\boldsymbol{\mathcal{B}}^{(2)}}^{\mathbf{f}^{(2)}}\rrbracket}
=
\left\llbracket
\mathrm{CH}^{(2)}_{\boldsymbol{\mathcal{B}}^{(2)}, \varphi(s)}\left(
f^{(2)\flat}_{s}\left(
\mathfrak{p}^{\mathbf{Pth}_{\boldsymbol{\mathcal{A}}^{(2)}}}
\right)
\right)
\right\rrbracket_{\varphi(s)}.
$$

\textsf{(4)}
For every $s\in S$, the interpretations of the operations $\mathrm{sc}_{s}^{0}$ and $\mathrm{tg}_{s}^{0}$ are given by
\begin{align*}
\mathrm{sc}_{s}^{0\llbracket\mathbf{PT}_{\boldsymbol{\mathcal{B}}^{(2)}}^{\mathbf{f}^{(2)}}\rrbracket}
&=
\mathrm{sc}_{\varphi(s)}^{0\llbracket\mathbf{PT}_{\boldsymbol{\mathcal{B}}^{(2)}}\rrbracket}
&&\mbox{and}&
\mathrm{tg}_{s}^{0\llbracket\mathbf{PT}_{\boldsymbol{\mathcal{B}}^{(2)}}^{\mathbf{f}^{(2)}}\rrbracket}
&=
\mathrm{tg}_{\varphi(s)}^{0\llbracket\mathbf{PT}_{\boldsymbol{\mathcal{B}}^{(2)}}\rrbracket}.
\end{align*}
That is, their interpretations are given by the interpretations of $\mathrm{sc}_{\varphi(s)}^{0}$ and $\mathrm{tg}_{\varphi(s)}^{0}$ in $\llbracket\mathbf{PT}_{\boldsymbol{\mathcal{B}}^{(2)}}\rrbracket$ that, we recall, was introduced in Proposition~\ref{PDPTQDCatAlg}.

\textsf{(5)}
Similarly,  for every $s \in S$, the interpretations of the operation $\circ_{s}^{0}$ is given by
$$
\circ_{s}^{0\llbracket\mathbf{PT}_{\boldsymbol{\mathcal{B}}^{(2)}}^{\mathbf{f}^{(2)}}\rrbracket}
=
\circ_{\varphi(s)}^{0\llbracket\mathbf{PT}_{\boldsymbol{\mathcal{B}}^{(2)}}\rrbracket}.
$$
That is, its interpretation is given by the interpretation of $\circ_{\varphi(s)}^{0}$ in $\llbracket\mathbf{PT}_{\boldsymbol{\mathcal{B}}^{(2)}}\rrbracket$ that, we recall, was introduced in Proposition~\ref{PDPTQDCatAlg}.

\textsf{(6)}
For every $s\in S$ and every $\mathfrak{p}^{(2)}\in\mathcal{A}_{s}^{(2)}$, the constant $\mathfrak{p}^{\mathbf{Pth}_{\boldsymbol{\mathcal{B}}^{(2)}}^{\mathbf{f}^{(2)}}}$ is given by
$$
\mathfrak{p}^{(2)\llbracket\mathbf{PT}_{\boldsymbol{\mathcal{B}}^{(2)}}^{\mathbf{f}^{(2)}}\rrbracket}
=
\left\llbracket
\mathrm{CH}^{(2)}_{\varphi(s)}\left(
f^{(2)}_{s}\left(
\mathfrak{p}^{(2)}
\right)
\right)
\right\rrbracket_{\varphi(s)}
\in
\left\llbracket\mathrm{PT}_{\boldsymbol{\mathcal{B}}^{(2)}}\right\rrbracket_{\varphi(s)}.
$$

\textsf{(7)}
For every $s\in S$, the interpretations of the operations $\mathrm{sc}_{s}^{1}$ and $\mathrm{tg}_{s}^{1}$ are given by
\allowdisplaybreaks
\begin{align*}
\mathrm{sc}_{s}^{1\llbracket\mathbf{PT}_{\boldsymbol{\mathcal{B}}^{(2)}}^{\mathbf{f}^{(2)}}\rrbracket}
&=
\mathrm{sc}_{\varphi(s)}^{1\llbracket\mathbf{PT}_{\boldsymbol{\mathcal{B}}^{(2)}}\rrbracket}
&&\mbox{and}&
\mathrm{tg}_{s}^{1\llbracket\mathbf{PT}_{\boldsymbol{\mathcal{B}}^{(2)}}^{\mathbf{f}^{(2)}}\rrbracket}
&=
\mathrm{tg}_{\varphi(s)}^{1\llbracket\mathbf{PT}_{\boldsymbol{\mathcal{B}}^{(2)}}\rrbracket}.
\end{align*}
That is, their interpretations are given by the interpretations of $\mathrm{sc}_{\varphi(s)}^{1}$ and $\mathrm{tg}_{\varphi(s)}^{1}$ in $\llbracket\mathbf{PT}_{\boldsymbol{\mathcal{B}}^{(2)}}\rrbracket$ that, we recall, was introduced in Proposition~\ref{PDPTQDCatAlg}.

\textsf{(8)}
Similarly, for every $s \in S$,  the interpretation of the partial binary operation $\circ_{s}^{1}$ is given by
$$
\circ_{s}^{1\llbracket\mathbf{PT}_{\boldsymbol{\mathcal{B}}^{(2)}}^{\mathbf{f}^{(2)}}\rrbracket}
=
\circ_{\varphi(s)}^{1\llbracket\mathbf{PT}_{\boldsymbol{\mathcal{B}}^{(2)}}\rrbracket}.
$$
That is, its interpretation is given by the interpretation of $\circ_{\varphi(s)}^{1}$ in $\llbracket\mathbf{PT}_{\boldsymbol{\mathcal{B}}^{(2)}}\rrbracket$ that, we recall, was introduced in Proposition~\ref{PDPTQDCatAlg}.

This completes the definition of the partial $\Sigma^{\boldsymbol{\mathcal{A}}^{(2)}}$-algebra $\llbracket\mathbf{PT}_{\boldsymbol{\mathcal{B}}^{(2)}}^{\mathbf{f}^{(2)}}\rrbracket$.
\end{proof}

\begin{proposition}
\label{PDPrPTBDCatHom}
Let $\mathbf{f}^{(2)}=(\varphi, c, (f^{(i)})_{i\in 3})$ be a second-order morphism from $\boldsymbol{\mathcal{A}}^{(2)}$ to $\boldsymbol{\mathcal{B}}^{(2)}$. Then the mapping 
$$
\mathrm{pr}_{\boldsymbol{\mathcal{B}}^{(2)}, \varphi}^{\Theta^{\llbracket 2 \rrbracket}}
\colon 
\mathrm{PT}_{\boldsymbol{\mathcal{B}}^{(2)}, \varphi}
\mor
\llbracket\mathrm{PT}_{\boldsymbol{\mathcal{B}}^{(2)}}\rrbracket_{\varphi}
$$
is a surjective $\Sigma^{\boldsymbol{\mathcal{A}}^{(2)}}$-homomorphism from $\mathbf{PT}_{\boldsymbol{\mathcal{B}}^{(2)}}^{\mathbf{f}^{(2)}}$ to $\llbracket\mathbf{PT}_{\boldsymbol{\mathcal{B}}^{(2)}}^{\mathbf{f}^{(2)}}\rrbracket$.
\end{proposition}

\begin{proof}
We prove that $\mathrm{pr}_{\boldsymbol{\mathcal{B}}^{(2)}, \varphi}^{\Theta^{\llbracket 2 \rrbracket}}$ is compatible with every operation symbol in $\Sigma^{\boldsymbol{\mathcal{A}}^{(2)}}$.

{\sffamily The mapping $\mathrm{pr}_{\boldsymbol{\mathcal{B}}^{(2)}, \varphi}^{\Theta^{\llbracket 2 \rrbracket}}$ is a $\Sigma$-homomorphism.}

Note that $\mathrm{pr}_{\boldsymbol{\mathcal{B}}^{(2)}, \varphi}^{\Theta^{\llbracket 2 \rrbracket}} = \mathbf{c}_{\mathfrak{d}}^{\ast} (\mathrm{pr}_{\boldsymbol{\mathcal{B}}^{(2)}}^{\Theta^{\llbracket 2 \rrbracket}})$. By Proposition~\ref{CDPTQPr}, the mapping $\mathrm{pr}_{\boldsymbol{\mathcal{B}}^{(2)}}^{\Theta^{\llbracket 2 \rrbracket}}$ is a $\Lambda^{\boldsymbol{\mathcal{B}}^{(2)}}$-homomorphism, thus in particular a $\Lambda$-homomorphism. Therefore, it follows from Proposition~\ref{PFunSig} that the mapping $\mathrm{pr}_{\boldsymbol{\mathcal{B}}^{(2)}, \varphi}^{\Theta^{\llbracket 2 \rrbracket}}$ is a $\Sigma$-homomorphism.


{\sffamily The mapping $\mathrm{pr}_{\boldsymbol{\mathcal{B}}^{(2)}, \varphi}^{\llbracket \cdot \rrbracket}$ is compatible with the first-order rewrite rules.}

Let $s$ be a sort in $S$ and $\mathfrak{p}$ a rewrite rule in $\mathcal{A}_{s}^{(1)}$. Thus,
$$
\left\llbracket
\mathfrak{p}^{\mathbf{PT}_{\boldsymbol{\mathcal{B}}^{(2)}}^{\mathbf{f}^{(2)}}}
\right\rrbracket_{\varphi(s)}
=
\left\llbracket
\mathrm{CH}^{(2)}_{\boldsymbol{\mathcal{B}}^{(2)}, \varphi(s)}\left(
f^{(2)\flat}_{s}\left(
\mathfrak{p}^{\mathbf{Pth}_{\boldsymbol{\mathcal{A}}^{(2)}}}
\right)
\right)
\right\rrbracket_{\varphi(s)}
=
\mathfrak{p}^{\llbracket\mathbf{PT}_{\boldsymbol{\mathcal{B}}^{(2)}}^{\mathbf{f}^{(2)}}\rrbracket}.
$$

Hence, $\mathrm{pr}_{\boldsymbol{\mathcal{B}}^{(2)}, \varphi}^{\Theta^{\llbracket 2 \rrbracket}_{\varphi(s)}}$ is compatible with the first-order rewrite rules.

{\sffamily The mapping $\mathrm{pr}_{\boldsymbol{\mathcal{B}}^{(2)}, \varphi}^{\Theta^{\llbracket 2 \rrbracket}}$ is compatible with the $0$-source.}

Let $s$ be a sort in $S$ and let us consider the $0$-source operation symbol $\mathrm{sc}_{s}^{0}$ in $\Sigma^{\boldsymbol{\mathcal{A}}^{(2)}}_{s,s}$. Let $P$ be a second-order path term in $\mathrm{PT}_{\boldsymbol{\mathcal{B}}^{(2)}, \varphi(s)}$.

The following chain of equalities holds
\allowdisplaybreaks
\begin{align*}
\left\llbracket
\mathrm{sc}_{s}^{0\mathbf{PT}_{\boldsymbol{\mathcal{B}}^{(2)}}^{\mathbf{f}^{(2)}}}\left(
P
\right)
\right\rrbracket_{\varphi(s)}
&=
\left\llbracket
\mathrm{sc}_{\varphi(s)}^{0\mathbf{PT}_{\boldsymbol{\mathcal{B}}^{(2)}}}\left(
P
\right)
\right\rrbracket_{\varphi(s)}
\tag{1}
\\
&=
\mathrm{sc}_{\varphi(s)}^{0\llbracket\mathbf{Pth}_{\boldsymbol{\mathcal{B}}^{(2)}}\rrbracket}\left(
\left\llbracket
P
\right\rrbracket_{\varphi(s)}
\right)
\tag{2}
\\
&=
\mathrm{sc}_{s}^{0\llbracket\mathbf{Pth}_{\boldsymbol{\mathcal{B}}^{(2)}}^{\mathbf{f}^{(2)}}\rrbracket}\left(
\left\llbracket
P
\right\rrbracket_{\varphi(s)}
\right).
\tag{3}
\end{align*}
The first equality unravels the interpretation of the operation symbol $\mathrm{sc}_{s}^{0}$ in the partial $\Sigma^{\boldsymbol{\mathcal{A}}^{(2)}}$-algebra $\mathbf{PT}_{\boldsymbol{\mathcal{B}}^{(2)}}^{\mathbf{f}^{(2)}}$, introduced in Proposition~\ref{PDPTBDCatAlg};
the second equality follows from the fact that, according to Proposition~\ref{CDPTQPr}, $\mathrm{pr}^{\Theta^{\llbracket 2 \rrbracket}}_{\boldsymbol{\mathcal{B}}^{(2)}}$ is a $\Lambda^{\boldsymbol{\mathcal{B}}^{(2)}}$-homomorphism from $\mathbf{PT}_{\boldsymbol{\mathcal{B}}^{(2)}}$ to $\llbracket\mathbf{PT}_{\boldsymbol{\mathcal{B}}^{(2)}}\rrbracket$;
finally, the last equality recovers the interpretation of the operation symbol $\mathrm{sc}_{s}^{0}$ in the partial $\Sigma^{\boldsymbol{\mathcal{A}}^{(2)}}$-algebra $\llbracket\mathbf{PT}_{\boldsymbol{\mathcal{B}}^{(2)}}^{\mathbf{f}^{(2)}}\rrbracket$, introduced in Proposition~\ref{PDQPTBDCatAlg}.

Hence, $\mathrm{pr}_{\boldsymbol{\mathcal{B}}^{(2)}, \varphi}^{\Theta^{\llbracket 2 \rrbracket}}$ is compatible with the $0$-source operation.

{\sffamily The mapping $\mathrm{pr}_{\boldsymbol{\mathcal{B}}^{(2)}, \varphi}^{\Theta^{\llbracket 2 \rrbracket}}$ is compatible with the $0$-target.}

Let $s$ be a sort in $S$ and let us consider the $0$-source operation symbol $\mathrm{tg}_{s}^{0}$ in $\Sigma^{\boldsymbol{\mathcal{A}}^{(2)}}_{s,s}$. Let $P$ be a second-order path term in $\mathrm{PT}_{\boldsymbol{\mathcal{B}}^{(2)}, \varphi(s)}$, then the following equality holds
$$
\left\llbracket
\mathrm{tg}_{s}^{0\mathbf{PT}_{\boldsymbol{\mathcal{B}}^{(2)}}^{\mathbf{f}^{(2)}}}\left(
P
\right)
\right\rrbracket_{\varphi(s)}
=
\mathrm{tg}_{s}^{0\llbracket\mathbf{PT}_{\boldsymbol{\mathcal{B}}^{(2)}}^{\mathbf{f}^{(2)}}\rrbracket}\left(
\left\llbracket
P
\right\rrbracket_{\varphi(s)}
\right).
$$

The proof of this case is identical to that of the $0$-source.

Hence, $\mathrm{pr}_{\boldsymbol{\mathcal{B}}^{(2)}, \varphi}^{\Theta^{\llbracket 2 \rrbracket}}$ is compatible with the $0$-target operation.

{\sffamily The mapping $\mathrm{pr}_{\boldsymbol{\mathcal{B}}^{(2)}, \varphi}^{\Theta^{\llbracket 2 \rrbracket}}$ is compatible with the $0$-composition.}

Let $s$ be a sort in $S$ and let us consider the $0$-composition operation symbol $\circ_{s}^{0}$ in $\Sigma_{ss,s}^{\boldsymbol{\mathcal{A}}^{(2)}}$. Let $P$ and $Q$ be two second-order path terms in $\mathrm{PT}_{\boldsymbol{\mathcal{B}}^{(2)}, \varphi(s)}$ such that
$$
\mathrm{sc}_{\boldsymbol{\mathcal{B}}^{(2)}, \varphi(s)}^{(1,2)}\left(
\mathrm{ip}^{(2,Y)@}_{\boldsymbol{\mathcal{B}}^{(2)}, \varphi(s)}\left(
Q
\right)
\right)
=
\mathrm{tg}_{\boldsymbol{\mathcal{B}}^{(2)}, \varphi(s)}^{(1,2)}\left(
\mathrm{ip}^{(2,Y)@}_{\boldsymbol{\mathcal{B}}^{(2)}, \varphi(s)}\left(
P
\right)
\right).
$$
Then the following equality holds
$$
\left\llbracket
Q
\circ_{s}^{1\mathbf{Pth}_{\boldsymbol{\mathcal{B}}^{(2)}}^{\mathbf{f}^{(2)}}}
P
\right\rrbracket_{\varphi(s)}
=
\left\llbracket
Q
\right\rrbracket_{\varphi(s)}
\circ_{s}^{1\llbracket\mathbf{Pth}_{\boldsymbol{\mathcal{B}}^{(2)}}^{\mathbf{f}^{(2)}}\rrbracket}
\left\llbracket
P
\right\rrbracket_{\varphi(s)}.
$$

The proof of this case is identical to that of the $0$-source.

Hence, $\mathrm{pr}_{\boldsymbol{\mathcal{B}}^{(2)}, \varphi}^{\Theta^{\llbracket 2 \rrbracket}}$ is compatible with the $0$-composition operation.

{\sffamily The mapping $\mathrm{pr}_{\boldsymbol{\mathcal{B}}^{(2)}, \varphi}^{\Theta^{\llbracket 2 \rrbracket}}$ is compatible with the second-order rewrite rules.}

Let $s$ be a sort in $S$ and $\mathfrak{p}^{(2)}$ a rewrite rule in $\mathcal{A}_{s}^{(2)}$. Thus,
$$
\left\llbracket
\mathfrak{p}^{(2)\mathbf{PT}_{\boldsymbol{\mathcal{B}}^{(2)}}^{\mathbf{f}^{(2)}}}
\right\rrbracket_{\varphi(s)}
=
\left\llbracket
\mathrm{CH}^{(2)}_{\boldsymbol{\mathcal{B}}^{(2)}, \varphi(s)}\left(
f_{s}^{(2)\flat}\left(
\mathfrak{p}^{(2)}
\right)
\right)
\right\rrbracket_{\varphi(s)}
=
\mathfrak{p}^{(2)\llbracket\mathbf{PT}_{\boldsymbol{\mathcal{B}}^{(2)}}^{\mathbf{f}^{(2)}}\rrbracket}.
$$

Hence, $\mathrm{pr}_{\boldsymbol{\mathcal{B}}^{(2)}, \varphi}^{\Theta^{\llbracket 2 \rrbracket}}$ is compatible with the second-order rewrite rules.

{\sffamily The mapping $\mathrm{pr}_{\boldsymbol{\mathcal{B}}^{(2)}, \varphi}^{\Theta^{\llbracket 2 \rrbracket}}$ is compatible with the $1$-source.}

Let $s$ be a sort in $S$ and let us consider the $1$-source operation symbol $\mathrm{sc}_{s}^{1}$ in $\Sigma^{\boldsymbol{\mathcal{A}}^{(2)}}_{s,s}$. Let $P$ be a second-order path term in $\mathrm{PT}_{\boldsymbol{\mathcal{B}}^{(2)}, \varphi(s)}$.

The following chain of equalities holds
\allowdisplaybreaks
\begin{align*}
\left\llbracket
\mathrm{sc}_{s}^{1\mathbf{PT}_{\boldsymbol{\mathcal{B}}^{(2)}}^{\mathbf{f}^{(2)}}}\left(
P
\right)
\right\rrbracket_{\varphi(s)}
&=
\left\llbracket
\mathrm{sc}_{\varphi(s)}^{1\mathbf{PT}_{\boldsymbol{\mathcal{B}}^{(2)}}}\left(
P
\right)
\right\rrbracket_{\varphi(s)}
\tag{1}
\\
&=
\mathrm{sc}_{\varphi(s)}^{1\llbracket\mathbf{PT}_{\boldsymbol{\mathcal{B}}^{(2)}}\rrbracket}\left(
\left\llbracket
P
\right\rrbracket_{\varphi(s)}
\right)
\tag{2}
\\
&=
\mathrm{sc}_{s}^{1\llbracket\mathbf{PT}_{\boldsymbol{\mathcal{B}}^{(2)}}^{\mathbf{f}^{(2)}}\rrbracket}\left(
\left\llbracket
P
\right\rrbracket_{\varphi(s)}
\right).
\tag{3}
\end{align*}
The first equality unravels the interpretation of the operation symbol $\mathrm{sc}_{s}^{1}$ in the partial $\Sigma^{\boldsymbol{\mathcal{A}}^{(2)}}$-algebra $\mathbf{PT}_{\boldsymbol{\mathcal{B}}^{(2)}}^{\mathbf{f}^{(2)}}$, introduced in Proposition~\ref{PDPTBDCatAlg};
the second equality follows from the fact that, according to Proposition~\ref{CDPTQPr}, $\mathrm{pr}^{\Theta^{\llbracket 2 \rrbracket}}_{\boldsymbol{\mathcal{B}}^{(2)}}$ is a $\Lambda^{\boldsymbol{\mathcal{B}}^{(2)}}$-homomorphism from $\mathbf{PT}_{\boldsymbol{\mathcal{B}}^{(2)}}$ to $\llbracket\mathbf{PT}_{\boldsymbol{\mathcal{B}}^{(2)}}\rrbracket$;
finally, the last equality recovers the interpretation of the operation symbol $\mathrm{sc}_{s}^{1}$ in the partial $\Sigma^{\boldsymbol{\mathcal{A}}^{(2)}}$-algebra $\llbracket\mathbf{PT}_{\boldsymbol{\mathcal{B}}^{(2)}}^{\mathbf{f}^{(2)}}\rrbracket$, introduced in Proposition~\ref{PDQPTBDCatAlg}.

Hence, $\mathrm{pr}_{\boldsymbol{\mathcal{B}}^{(2)}, \varphi}^{\Theta^{\llbracket 2 \rrbracket}}$ is compatible with the $1$-source operation.

{\sffamily The mapping $\mathrm{pr}_{\boldsymbol{\mathcal{B}}^{(2)}, \varphi}^{\Theta^{\llbracket 2 \rrbracket}}$ is compatible with the $1$-target.}

Let $s$ be a sort in $S$ and let us consider the $1$-source operation symbol $\mathrm{sc}_{s}^{1}$ in $\Sigma^{\boldsymbol{\mathcal{A}}^{(2)}}_{s,s}$. Let $P$ be a second-order path term in $\mathrm{PT}_{\boldsymbol{\mathcal{B}}^{(2)}, \varphi(s)}$, then the following equality holds
$$
\left\llbracket
\mathrm{tg}_{s}^{1\mathbf{PT}_{\boldsymbol{\mathcal{B}}^{(2)}}^{\mathbf{f}^{(2)}}}\left(
P
\right)
\right\rrbracket_{\varphi(s)}
=
\mathrm{tg}_{s}^{1\llbracket\mathbf{PT}_{\boldsymbol{\mathcal{B}}^{(2)}}^{\mathbf{f}^{(2)}}\rrbracket}\left(
\left\llbracket
P
\right\rrbracket_{\varphi(s)}
\right).
$$

The proof of this case is identical to that of the $1$-source.

Hence, $\mathrm{pr}_{\boldsymbol{\mathcal{B}}^{(2)}, \varphi}^{\Theta^{\llbracket 2 \rrbracket}}$ is compatible with the $1$-target operation.

{\sffamily The mapping $\mathrm{pr}_{\boldsymbol{\mathcal{B}}^{(2)}, \varphi}^{\Theta^{\llbracket 2 \rrbracket}}$ is compatible with the $1$-composition.}

Let $s$ be a sort in $S$ and let us consider the $1$-composition operation symbol $\circ_{s}^{1}$ in $\Sigma_{ss,s}^{\boldsymbol{\mathcal{A}}^{(2)}}$. Let $P$ and $Q$ be two second-order path terms in $\mathrm{PT}_{\boldsymbol{\mathcal{B}}^{(2)}, \varphi(s)}$ such that
$$
\mathrm{sc}_{\boldsymbol{\mathcal{B}}^{(2)}, \varphi(s)}^{([1],2)}\left(
\mathrm{ip}^{(2,Y)@}_{\boldsymbol{\mathcal{B}}^{(2)}, \varphi(s)}\left(
Q
\right)
\right)
=
\mathrm{tg}_{\boldsymbol{\mathcal{B}}^{(2)}, \varphi(s)}^{([1],2)}\left(
\mathrm{ip}^{(2,Y)@}_{\boldsymbol{\mathcal{B}}^{(2)}, \varphi(s)}\left(
P
\right)
\right).
$$
Then the following equality holds
$$
\left\llbracket
Q
\circ_{s}^{1\mathbf{PT}_{\boldsymbol{\mathcal{B}}^{(2)}}^{\mathbf{f}^{(2)}}}
P
\right\rrbracket_{\varphi(s)}
=
\left\llbracket
Q
\right\rrbracket_{\varphi(s)}
\circ_{s}^{1\llbracket\mathbf{PT}_{\boldsymbol{\mathcal{B}}^{(2)}}^{\mathbf{f}^{(2)}}\rrbracket}
\left\llbracket
P
\right\rrbracket_{\varphi(s)}.
$$

The proof of this case is identical to that of the $1$-source.

Hence, $\mathrm{pr}_{\boldsymbol{\mathcal{B}}^{(2)}, \varphi}^{\Theta^{\llbracket 2 \rrbracket}}$ is compatible with the $1$-composition operation.

This completes the proof.
\end{proof}

Finally, we show that the partial $\Sigma^{\boldsymbol{\mathcal{A}}^{(2)}}$-algebras $\llbracket\mathbf{Pth}_{\boldsymbol{\mathcal{B}}^{(2)}}^{\mathbf{f}^{(2)}}\rrbracket$ and $\llbracket\mathbf{PT}_{\boldsymbol{\mathcal{B}}^{(2)}}^{\mathbf{f}^{(2)}}\rrbracket$ are isomorphic. Indeed the mappings $\mathrm{ip}^{(\llbracket 2 \rrbracket, Y)@}_{\boldsymbol{\mathcal{B}}^{(2)}, \varphi}$ and $\mathrm{CH}^{\llbracket 2 \rrbracket}_{\boldsymbol{\mathcal{B}}^{(2)}, \varphi}$ are a pair of inverse $\Sigma^{\boldsymbol{\mathcal{A}}^{(2)}}$-isomorphisms.

\begin{proposition}
\label{PDQCHBDCatHom}
Let $\mathbf{f}^{(2)}=(\varphi, c, (f^{(i)})_{i\in 3})$ be a second-order morphism from $\boldsymbol{\mathcal{A}}^{(2)}$ to $\boldsymbol{\mathcal{B}}^{(2)}$. Then the mapping 
$$
\mathrm{CH}^{\llbracket2\rrbracket}_{\boldsymbol{\mathcal{B}}^{(2)}, \varphi}
\colon 
\llbracket\mathrm{Pth}_{\boldsymbol{\mathcal{B}}^{(2)}}\rrbracket_{\varphi}
\mor
\llbracket\mathrm{PT}_{\boldsymbol{\mathcal{B}}^{(2)}}\rrbracket_{\varphi}
$$
is a $\Sigma^{\boldsymbol{\mathcal{A}}^{(2)}}$-homomorphism from $\llbracket\mathbf{Pth}_{\boldsymbol{\mathcal{B}}^{(2)}}^{\mathbf{f}^{(2)}}\rrbracket$ to $\llbracket\mathbf{PT}_{\boldsymbol{\mathcal{B}}^{(2)}}^{\mathbf{f}^{(2)}}\rrbracket$.
\end{proposition}

\begin{proof}
We prove that $\mathrm{CH}^{\llbracket2\rrbracket}_{\boldsymbol{\mathcal{B}}^{(2)}, \varphi}$ is compatible with every operation symbol in $\Sigma^{\boldsymbol{\mathcal{A}}^{(2)}}$.

{\sffamily The mapping $\mathrm{CH}^{\llbracket2\rrbracket}_{\boldsymbol{\mathcal{B}}^{(2)}, \varphi}$ is a $\Sigma$-homomorphism.}

Note that $\mathrm{CH}^{\llbracket2\rrbracket}_{\boldsymbol{\mathcal{B}}^{(2)}, \varphi} = \mathbf{c}_{\mathfrak{d}}^{\ast} (\mathrm{CH}^{\llbracket2\rrbracket}_{\boldsymbol{\mathcal{B}}^{(2)}})$. By Proposition~\ref{CDIso}, the mapping $\mathrm{CH}^{\llbracket2\rrbracket}_{\boldsymbol{\mathcal{B}}^{(2)}}$ is a $\Lambda^{\boldsymbol{\mathcal{B}}^{(2)}}$-homomorphism, thus in particular a $\Lambda$-homomorphism. Therefore, it follows from Proposition~\ref{PFunSig} that the mapping $\mathrm{CH}^{\llbracket2\rrbracket}_{\boldsymbol{\mathcal{B}}^{(2)}, \varphi}$ is a $\Sigma$-homomorphism.




{\sffamily The mapping $\mathrm{CH}^{\llbracket2\rrbracket}_{\boldsymbol{\mathcal{B}}^{(2)}, \varphi}$ is compatible with the first-order rewrite rules.}

Let $s$ be a sort in $S$ and $\mathfrak{p}$ a rewrite rule in $\mathcal{A}_{s}^{(1)}$. Thus,
\begin{align*}
\mathrm{CH}^{\llbracket2\rrbracket}_{\boldsymbol{\mathcal{B}}^{(2)}, \varphi(s)}\left(
\mathfrak{p}^{\llbracket\mathbf{Pth}_{\boldsymbol{\mathcal{B}}^{(2)}}^{\mathbf{f}^{(2)}}\rrbracket}
\right)
&=
\mathrm{CH}^{\llbracket2\rrbracket}_{\boldsymbol{\mathcal{B}}^{(2)}, \varphi(s)}\left(
\left\llbracket
f^{(2)\flat}\left(
\mathrm{ech}^{(2,\mathcal{A}^{(1)})}_{\boldsymbol{\mathcal{A}}^{(2)}, s} \left(
\mathfrak{p}
\right)
\right)
\right\rrbracket_{\varphi(s)}
\right)
\tag{1}
\\
&=
\left\llbracket
\mathrm{CH}^{(2)}_{\boldsymbol{\mathcal{B}}^{(2)}, \varphi(s)}\left(
f^{(2)\flat}\left(
\mathrm{ech}^{(2,\mathcal{A}^{(1)})}_{\boldsymbol{\mathcal{A}}^{(2)}, s} \left(
\mathfrak{p}
\right)
\right)
\right)
\right\rrbracket_{\varphi(s)}
\tag{2}
\\
&=
\mathfrak{p}^{\llbracket\mathbf{PT}_{\boldsymbol{\mathcal{B}}^{(2)}}^{\mathbf{f}^{(2)}}\rrbracket}
\tag{3}
\end{align*}
The first equality unravels the definition of the constant $\mathfrak{p}$ in the partial $\Sigma^{\boldsymbol{\mathcal{A}}^{(2)}}$-algebra $\llbracket\mathbf{Pth}_{\boldsymbol{\mathcal{B}}^{(2)}}^{\mathbf{f}^{(2)}}\rrbracket$, introduced in Proposition~\ref{PDQPthBDCatAlg};
the second equality unravels the definition of the mapping $\mathrm{CH}^{\llbracket2\rrbracket}_{\boldsymbol{\mathcal{B}}^{(2)}}$, introduced in Definition~\ref{DDPTQDCH};
finally, the last equality recovers the definition of the constant $\mathfrak{p}$ in the partial $\Sigma^{\boldsymbol{\mathcal{A}}^{(2)}}$-algebra $\llbracket\mathbf{PT}_{\boldsymbol{\mathcal{B}}^{(2)}}^{\mathbf{f}^{(2)}}\rrbracket$, introduced in Proposition~\ref{PDQPTBDCatAlg};

Hence, $\mathrm{CH}^{\llbracket2\rrbracket}_{\boldsymbol{\mathcal{B}}^{(2)}, \varphi}$ is compatible with the first-order rewrite rules.

{\sffamily The mapping $\mathrm{CH}^{\llbracket2\rrbracket}_{\boldsymbol{\mathcal{B}}^{(2)}, \varphi}$ is compatible with the $0$-source.}

Let $s$ be a sort in $S$ and let us consider the $0$-source operation symbol $\mathrm{sc}_{s}^{0}$ in $\Sigma^{\boldsymbol{\mathcal{A}}^{(1)}}_{s,s}$. Let $\mathfrak{P}^{(2)}$ be a second-order path in $\mathrm{Pth}_{\boldsymbol{\mathcal{B}}^{(2)}, \varphi(s)}$.

The following chain of equalities holds
\allowdisplaybreaks
\begin{flushleft}
$
\mathrm{CH}^{\llbracket2\rrbracket}_{\boldsymbol{\mathcal{B}}^{(2)}, \varphi(s)} \left(
\mathrm{sc}_{s}^{0\llbracket\mathbf{Pth}_{\boldsymbol{\mathcal{B}}^{(2)}}^{\mathbf{f}^{(2)}}\rrbracket}\left(
\left\llbracket
\mathfrak{P}^{(2)}
\right\rrbracket_{\varphi(s)}
\right)
\right)
$
\begin{align*}
&=
\mathrm{CH}^{\llbracket2\rrbracket}_{\boldsymbol{\mathcal{B}}^{(2)}, \varphi(s)} \left(
\mathrm{sc}_{\varphi(s)}^{0\llbracket\mathbf{Pth}_{\boldsymbol{\mathcal{B}}^{(2)}}\rrbracket}\left(
\left\llbracket
\mathfrak{P}^{(2)}
\right\rrbracket_{\varphi(s)}
\right)
\right)
\tag{1}
\\
&=
\mathrm{sc}_{\varphi(s)}^{0\llbracket\mathbf{PT}_{\boldsymbol{\mathcal{B}}^{(2)}}\rrbracket}\left(
\mathrm{CH}^{\llbracket2\rrbracket}_{\boldsymbol{\mathcal{B}}^{(2)}, \varphi(s)} \left(
\left\llbracket
\mathfrak{P}^{(2)}
\right\rrbracket_{\varphi(s)}
\right)
\right)
\tag{2}
\\
&=
\mathrm{sc}_{s}^{0\llbracket\mathbf{PT}_{\boldsymbol{\mathcal{B}}^{(2)}}^{\mathbf{f}^{(2)}}\rrbracket}\left(
\mathrm{CH}^{\llbracket2\rrbracket}_{\boldsymbol{\mathcal{B}}^{(2)}, \varphi(s)} \left(
\left\llbracket
\mathfrak{P}^{(2)}
\right\rrbracket_{\varphi(s)}
\right)
\right).
\tag{3}
\end{align*}
\end{flushleft}

The first equality unravels the interpretation of the operation symbol $\mathrm{sc}_{s}^{0}$ in the partial $\Sigma^{\boldsymbol{\mathcal{A}}^{(1)}}$-algebra $\llbracket\mathbf{Pth}_{\boldsymbol{\mathcal{B}}^{(2)}}^{\mathbf{f}^{(2)}}\rrbracket$, introduced in Proposition~\ref{PDPthBDCatAlg};
the second equality follows from the fact that, according to Claim~\ref{CDIso}, $\mathrm{CH}^{\llbracket2\rrbracket}_{\boldsymbol{\mathcal{B}}^{(2)}, \varphi}$ is a $\Lambda^{\boldsymbol{\mathcal{B}}^{(2)}}$-homomorphism from $\llbracket\mathbf{Pth}_{\boldsymbol{\mathcal{B}}^{(2)}}\rrbracket$ to $\llbracket\mathbf{PT}_{\boldsymbol{\mathcal{B}}^{(2)}}\rrbracket$;
finally, the last equality recovers the interpretation of the operation symbol $\mathrm{sc}_{s}^{0}$ in the partial $\Sigma^{\boldsymbol{\mathcal{A}}^{(2)}}$-algebra $\llbracket\mathbf{PT}_{\boldsymbol{\mathcal{B}}^{(2)}}^{\mathbf{f}^{(2)}}\rrbracket$, introduced in Proposition~\ref{PDQPTBDCatAlg}.

Hence, $\mathrm{CH}^{\llbracket2\rrbracket}_{\boldsymbol{\mathcal{B}}^{(2)}, \varphi}$ is compatible with the $0$-source operation.

{\sffamily The mapping $\mathrm{CH}^{\llbracket2\rrbracket}_{\boldsymbol{\mathcal{B}}^{(2)}, \varphi}$ is compatible with the $0$-target.}

Let $s$ be a sort in $S$ and let us consider the $0$-target operation symbol $\mathrm{tg}_{s}^{0}$ in $\Sigma^{\boldsymbol{\mathcal{A}}^{(1)}}_{s,s}$. Let $\mathfrak{P}^{(2)}$ be a second-order path in $\mathrm{Pth}_{\boldsymbol{\mathcal{B}}^{(2)}, \varphi(s)}$, then the following equality holds
\begin{multline*}
\mathrm{CH}^{\llbracket2\rrbracket}_{\boldsymbol{\mathcal{B}}^{(2)}, \varphi(s)}\left(
\mathrm{tg}_{s}^{0\llbracket\mathbf{Pth}_{\boldsymbol{\mathcal{B}}^{(2)}}^{\mathbf{f}^{(2)}}\rrbracket}\left(
\left\llbracket
\mathfrak{P}^{(2)}
\right\rrbracket_{\varphi(s)}
\right)
\right)
\\
=
\mathrm{tg}_{s}^{0\llbracket\mathbf{PT}_{\boldsymbol{\mathcal{B}}^{(2)}}^{\mathbf{f}^{(2)}}\rrbracket}\left(
\mathrm{CH}^{\llbracket2\rrbracket}_{\boldsymbol{\mathcal{B}}^{(2)}, \varphi(s)}\left(
\left\llbracket
\mathfrak{P}^{(2)}
\right\rrbracket_{\varphi(s)}
\right)
\right).
\end{multline*}

The proof of this case is identical to that of the $0$-source.

Hence, $\mathrm{CH}^{\llbracket2\rrbracket}_{\boldsymbol{\mathcal{B}}^{(2)}, \varphi}$ is compatible with the $0$-target operation.

{\sffamily The mapping $\mathrm{CH}^{\llbracket2\rrbracket}_{\boldsymbol{\mathcal{B}}^{(2)}, \varphi}$ is compatible with the $0$-composition.}

Let $s$ be a sort in $S$ and let us consider the $0$-composition operation symbol $\circ_{s}^{0}$ in $\Sigma_{ss,s}^{\boldsymbol{\mathcal{A}}^{(1)}}$. Let $\mathfrak{P}^{(2)}$ and $\mathfrak{Q}^{(2)}$ be two second-order paths in $\mathrm{Pth}_{\boldsymbol{\mathcal{B}}^{(2)}, \varphi(s)}$ such that
$$
\mathrm{sc}_{\boldsymbol{\mathcal{B}}^{(2)}, \varphi(s)}^{(0,2)}\left(\mathfrak{Q}^{(2)}\right)
=
\mathrm{tg}_{\boldsymbol{\mathcal{B}}^{(2)}, \varphi(s)}^{(0,2)}\left(\mathfrak{P}^{(2)}\right).
$$

Then the following equality holds
\begin{multline*}
\mathrm{CH}^{\llbracket2\rrbracket}_{\boldsymbol{\mathcal{B}}^{(2)}, \varphi(s)}\left(
\left\llbracket
\mathfrak{Q}^{(2)}
\right\rrbracket_{\varphi(s)}
\circ_{s}^{0\llbracket\mathbf{Pth}_{\boldsymbol{\mathcal{B}}^{(2)}}^{\mathbf{f}^{(2)}}\rrbracket}
\left\llbracket
\mathfrak{P}^{(2)}
\right\rrbracket_{\varphi(s)}
\right)
\\
=
\mathrm{CH}^{\llbracket2\rrbracket}_{\boldsymbol{\mathcal{B}}^{(2)}, \varphi(s)}\left(
\left\llbracket
\mathfrak{Q}^{(2)}
\right\rrbracket_{\varphi(s)}
\right)
\circ_{s}^{0\llbracket\mathbf{PT}_{\boldsymbol{\mathcal{B}}^{(2)}}^{\mathbf{f}^{(2)}}\rrbracket}
\mathrm{CH}^{\llbracket2\rrbracket}_{\boldsymbol{\mathcal{B}}^{(2)}, \varphi(s)}\left(
\left\llbracket
\mathfrak{P}^{(2)}
\right\rrbracket_{\varphi(s)}
\right).
\end{multline*}

The proof of this case is identical to that of the $0$-source.

Hence, $\mathrm{CH}^{\llbracket2\rrbracket}_{\boldsymbol{\mathcal{B}}^{(2)}, \varphi}$ is compatible with the $0$-composition operation.

{\sffamily The mapping $\mathrm{CH}^{\llbracket2\rrbracket}_{\boldsymbol{\mathcal{B}}^{(2)}, \varphi}$ is compatible with the second-order rewrite rules.}

Let $s$ be a sort in $S$ and $\mathfrak{p}^{(2)}$ a second-order rewrite rule in $\mathcal{A}_{s}^{(2)}$. Thus,
\begin{align*}
\mathrm{CH}^{\llbracket2\rrbracket}_{\boldsymbol{\mathcal{B}}^{(2)}, \varphi(s)}\left(
\mathfrak{p}^{(2)\llbracket\mathbf{Pth}_{\boldsymbol{\mathcal{B}}^{(2)}}^{\mathbf{f}^{(2)}}\rrbracket}
\right)
&=
\mathrm{CH}^{\llbracket2\rrbracket}_{\boldsymbol{\mathcal{B}}^{(2)}, \varphi(s)}\left(
\left\llbracket
f^{(2)}\left(
\mathfrak{p}^{(2)}
\right)
\right\rrbracket_{\varphi(s)}
\right)
\tag{1}
\\
&=
\left\llbracket
\mathrm{CH}^{(2)}_{\boldsymbol{\mathcal{B}}^{(2)}, \varphi(s)}\left(
f^{(2)}\left(
\mathfrak{p}^{(2)}
\right)
\right)
\right\rrbracket_{\varphi(s)}
\tag{2}
\\
&=
\mathfrak{p}^{(2)\llbracket\mathbf{PT}_{\boldsymbol{\mathcal{B}}^{(2)}}^{\mathbf{f}^{(2)}}\rrbracket}
\tag{3}
\end{align*}
The first equality unravels the definition of the constant $\mathfrak{p}^{(2)}$ in the partial $\Sigma^{\boldsymbol{\mathcal{A}}^{(2)}}$-algebra $\llbracket\mathbf{Pth}_{\boldsymbol{\mathcal{B}}^{(2)}}^{\mathbf{f}^{(2)}}\rrbracket$, introduced in Proposition~\ref{PDQPthBDCatAlg};
the second equality unravels the definition of the mapping $\mathrm{CH}^{\llbracket2\rrbracket}_{\boldsymbol{\mathcal{B}}^{(2)}}$, introduced in Definition~\ref{DDPTQDCH};
finally, the last equality recovers the definition of the constant $\mathfrak{p}$ in the partial $\Sigma^{\boldsymbol{\mathcal{A}}^{(2)}}$-algebra $\llbracket\mathbf{PT}_{\boldsymbol{\mathcal{B}}^{(2)}}^{\mathbf{f}^{(2)}}\rrbracket$, introduced in Proposition~\ref{PDQPTBDCatAlg};

Hence, $\mathrm{CH}^{\llbracket2\rrbracket}_{\boldsymbol{\mathcal{B}}^{(2)}, \varphi}$ is compatible with the first-order rewrite rules.

{\sffamily The mapping $\mathrm{CH}^{\llbracket2\rrbracket}_{\boldsymbol{\mathcal{B}}^{(2)}, \varphi}$ is compatible with the $1$-source.}

Let $s$ be a sort in $S$ and let us consider the $1$-source operation symbol $\mathrm{sc}_{s}^{1}$ in $\Sigma^{\boldsymbol{\mathcal{A}}^{(1)}}_{s,s}$. Let $\mathfrak{P}^{(2)}$ be a second-order path in $\mathrm{Pth}_{\boldsymbol{\mathcal{B}}^{(2)}, \varphi(s)}$.

The following chain of equalities holds
\allowdisplaybreaks
\begin{flushleft}
$
\mathrm{CH}^{\llbracket2\rrbracket}_{\boldsymbol{\mathcal{B}}^{(2)}, \varphi(s)} \left(
\mathrm{sc}_{s}^{1\llbracket\mathbf{Pth}_{\boldsymbol{\mathcal{B}}^{(2)}}^{\mathbf{f}^{(2)}}\rrbracket}\left(
\left\llbracket
\mathfrak{P}^{(2)}
\right\rrbracket_{\varphi(s)}
\right)
\right)
$
\begin{align*}
&=
\mathrm{CH}^{\llbracket2\rrbracket}_{\boldsymbol{\mathcal{B}}^{(2)}, \varphi(s)} \left(
\mathrm{sc}_{\varphi(s)}^{1\llbracket\mathbf{Pth}_{\boldsymbol{\mathcal{B}}^{(2)}}\rrbracket}\left(
\left\llbracket
\mathfrak{P}^{(2)}
\right\rrbracket_{\varphi(s)}
\right)
\right)
\tag{1}
\\
&=
\mathrm{sc}_{\varphi(s)}^{1\llbracket\mathbf{PT}_{\boldsymbol{\mathcal{B}}^{(2)}}\rrbracket}\left(
\mathrm{CH}^{\llbracket2\rrbracket}_{\boldsymbol{\mathcal{B}}^{(2)}, \varphi(s)} \left(
\left\llbracket
\mathfrak{P}^{(2)}
\right\rrbracket_{\varphi(s)}
\right)
\right)
\tag{2}
\\
&=
\mathrm{sc}_{s}^{1\llbracket\mathbf{PT}_{\boldsymbol{\mathcal{B}}^{(2)}}^{\mathbf{f}^{(2)}}\rrbracket}\left(
\mathrm{CH}^{\llbracket2\rrbracket}_{\boldsymbol{\mathcal{B}}^{(2)}, \varphi(s)} \left(
\left\llbracket
\mathfrak{P}^{(2)}
\right\rrbracket_{\varphi(s)}
\right)
\right).
\tag{3}
\end{align*}
\end{flushleft}

The first equality unravels the interpretation of the operation symbol $\mathrm{sc}_{s}^{1}$ in the partial $\Sigma^{\boldsymbol{\mathcal{A}}^{(1)}}$-algebra $\llbracket\mathbf{Pth}_{\boldsymbol{\mathcal{B}}^{(2)}}^{\mathbf{f}^{(2)}}\rrbracket$, introduced in Proposition~\ref{PDPthBDCatAlg};
the second equality follows from the fact that, according to Claim~\ref{CDIso}, $\mathrm{CH}^{\llbracket2\rrbracket}_{\boldsymbol{\mathcal{B}}^{(2)}, \varphi}$ is a $\Lambda^{\boldsymbol{\mathcal{B}}^{(2)}}$-homomorphism from $\llbracket\mathbf{Pth}_{\boldsymbol{\mathcal{B}}^{(2)}}\rrbracket$ to $\llbracket\mathbf{PT}_{\boldsymbol{\mathcal{B}}^{(2)}}\rrbracket$;
finally, the last equality recovers the interpretation of the operation symbol $\mathrm{sc}_{s}^{1}$ in the partial $\Sigma^{\boldsymbol{\mathcal{A}}^{(2)}}$-algebra $\llbracket\mathbf{PT}_{\boldsymbol{\mathcal{B}}^{(2)}}^{\mathbf{f}^{(2)}}\rrbracket$, introduced in Proposition~\ref{PDQPTBDCatAlg}.

Hence, $\mathrm{CH}^{\llbracket2\rrbracket}_{\boldsymbol{\mathcal{B}}^{(2)}, \varphi}$ is compatible with the $1$-source operation.

{\sffamily The mapping $\mathrm{CH}^{\llbracket2\rrbracket}_{\boldsymbol{\mathcal{B}}^{(2)}, \varphi}$ is compatible with the $1$-target.}

Let $s$ be a sort in $S$ and let us consider the $1$-target operation symbol $\mathrm{tg}_{s}^{1}$ in $\Sigma^{\boldsymbol{\mathcal{A}}^{(1)}}_{s,s}$. Let $\mathfrak{P}^{(2)}$ be a second-order path in $\mathrm{Pth}_{\boldsymbol{\mathcal{B}}^{(2)}, \varphi(s)}$, then the following equality holds
\begin{multline*}
\mathrm{CH}^{\llbracket2\rrbracket}_{\boldsymbol{\mathcal{B}}^{(2)}, \varphi(s)}\left(
\mathrm{tg}_{s}^{1\llbracket\mathbf{Pth}_{\boldsymbol{\mathcal{B}}^{(2)}}^{\mathbf{f}^{(2)}}\rrbracket}\left(
\left\llbracket
\mathfrak{P}^{(2)}
\right\rrbracket_{\varphi(s)}
\right)
\right)
\\
=
\mathrm{tg}_{s}^{1\llbracket\mathbf{PT}_{\boldsymbol{\mathcal{B}}^{(2)}}^{\mathbf{f}^{(2)}}\rrbracket}\left(
\mathrm{CH}^{\llbracket2\rrbracket}_{\boldsymbol{\mathcal{B}}^{(2)}, \varphi(s)}\left(
\left\llbracket
\mathfrak{P}^{(2)}
\right\rrbracket_{\varphi(s)}
\right)
\right).
\end{multline*}

The proof of this case is identical to that of the $1$-source.

Hence, $\mathrm{CH}^{\llbracket2\rrbracket}_{\boldsymbol{\mathcal{B}}^{(2)}, \varphi}$ is compatible with the $1$-target operation.

{\sffamily The mapping $\mathrm{CH}^{\llbracket2\rrbracket}_{\boldsymbol{\mathcal{B}}^{(2)}, \varphi}$ is compatible with the $1$-composition.}

Let $s$ be a sort in $S$ and let us consider the $1$-composition operation symbol $\circ_{s}^{1}$ in $\Sigma_{ss,s}^{\boldsymbol{\mathcal{A}}^{(1)}}$. Let $\mathfrak{P}^{(2)}$ and $\mathfrak{Q}^{(2)}$ be two second-order paths in $\mathrm{Pth}_{\boldsymbol{\mathcal{B}}^{(2)}, \varphi(s)}$ such that
$$
\mathrm{sc}_{\boldsymbol{\mathcal{B}}^{(2)}, \varphi(s)}^{([1],2)}\left(\mathfrak{Q}^{(2)}\right)
=
\mathrm{tg}_{\boldsymbol{\mathcal{B}}^{(2)}, \varphi(s)}^{([1],2)}\left(\mathfrak{P}^{(2)}\right).
$$

Then the following equality holds
\begin{multline*}
\mathrm{CH}^{\llbracket2\rrbracket}_{\boldsymbol{\mathcal{B}}^{(2)}, \varphi(s)}\left(
\left\llbracket
\mathfrak{Q}^{(2)}
\right\rrbracket_{\varphi(s)}
\circ_{s}^{1\llbracket\mathbf{Pth}_{\boldsymbol{\mathcal{B}}^{(2)}}^{\mathbf{f}^{(2)}}\rrbracket}
\left\llbracket
\mathfrak{P}^{(2)}
\right\rrbracket_{\varphi(s)}
\right)
\\
=
\mathrm{CH}^{\llbracket2\rrbracket}_{\boldsymbol{\mathcal{B}}^{(2)}, \varphi(s)}\left(
\left\llbracket
\mathfrak{Q}^{(2)}
\right\rrbracket_{\varphi(s)}
\right)
\circ_{s}^{1\llbracket\mathbf{PT}_{\boldsymbol{\mathcal{B}}^{(2)}}^{\mathbf{f}^{(2)}}\rrbracket}
\mathrm{CH}^{\llbracket2\rrbracket}_{\boldsymbol{\mathcal{B}}^{(2)}, \varphi(s)}\left(
\left\llbracket
\mathfrak{P}^{(2)}
\right\rrbracket_{\varphi(s)}
\right).
\end{multline*}

The proof of this case is identical to that of the $1$-source.

Hence, $\mathrm{CH}^{\llbracket2\rrbracket}_{\boldsymbol{\mathcal{B}}^{(2)}, \varphi}$ is compatible with the $1$-composition operation.

This completes the proof.
\end{proof}

\begin{proposition}
\label{PDQIpfcBDCatHom}
Let $\mathbf{f}^{(2)}=(\varphi, c, (f^{(i)})_{i\in 3})$ be a second-order morphism from $\boldsymbol{\mathcal{A}}^{(1)}$ to $\boldsymbol{\mathcal{B}}^{(1)}$. Then the mapping 
$$
\mathrm{ip}^{(\llbracket2\rrbracket, Y)@}_{\boldsymbol{\mathcal{B}}^{(2)}, \varphi}
\colon 
\llbracket\mathrm{PT}_{\boldsymbol{\mathcal{B}}^{(2)}}\rrbracket_{\varphi}
\mor
\llbracket\mathrm{Pth}_{\boldsymbol{\mathcal{B}}^{(2)}}\rrbracket_{\varphi}
$$
is a $\Sigma^{\boldsymbol{\mathcal{A}}^{(2)}}$-homomorphism from $\llbracket\mathbf{PT}_{\boldsymbol{\mathcal{B}}^{(2)}}^{\mathbf{f}^{(2)}}\rrbracket$ to $\llbracket\mathbf{Pth}_{\boldsymbol{\mathcal{B}}^{(2)}}^{\mathbf{f}^{(2)}}\rrbracket$.
\end{proposition}

\begin{proof}
We prove that $\mathrm{ip}^{(\llbracket2\rrbracket, Y)@}_{\boldsymbol{\mathcal{B}}^{(2)}, \varphi}$ is compatible with every operation symbol in $\Sigma^{\boldsymbol{\mathcal{A}}^{(2)}}$.

{\sffamily The mapping $\mathrm{ip}^{(\llbracket2\rrbracket, Y)@}_{\boldsymbol{\mathcal{B}}^{(2)}, \varphi}$ is a $\Sigma$-homomorphism.}

Note that $\mathrm{ip}^{(\llbracket2\rrbracket, Y)@}_{\boldsymbol{\mathcal{B}}^{(2)}, \varphi} = \mathbf{c}_{\mathfrak{d}}^{\ast} (\mathrm{ip}^{(\llbracket2\rrbracket, Y)@}_{\boldsymbol{\mathcal{B}}^{(2)}})$. By Proposition~\ref{CDIsoIpfc}, the mapping $\mathrm{ip}^{(\llbracket2\rrbracket, Y)@}_{\boldsymbol{\mathcal{B}}^{(2)}}$ is a $\Lambda^{\boldsymbol{\mathcal{B}}^{(2)}}$-homomorphism, thus in particular a $\Lambda$-homomorphism. Therefore, it follows from Proposition~\ref{PFunSig} that the mapping $\mathrm{ip}^{(\llbracket2\rrbracket, Y)@}_{\boldsymbol{\mathcal{B}}^{(2)}, \varphi}$ is a $\Sigma$-homomorphism.


{\sffamily The mapping $\mathrm{ip}^{(\llbracket2\rrbracket, X)@}_{\boldsymbol{\mathcal{B}}^{(2)}, \varphi}$ is compatible with the first-order rewrite rules.}

Let $s$ be a sort in $S$ and $\mathfrak{p}$ a rewrite rule in $\mathcal{A}_{s}^{(1)}$. Thus,
\begin{flushleft}
$
\mathrm{ip}^{(\llbracket2\rrbracket, Y)@}_{\boldsymbol{\mathcal{B}}^{(2)}, \varphi(s)}\left(
\mathfrak{p}^{\llbracket\mathbf{PT}_{\boldsymbol{\mathcal{B}}^{(2)}}^{\mathbf{f}^{(2)}}\rrbracket}
\right)
$
\begin{align*}
&=
\mathrm{ip}^{(\llbracket2\rrbracket, Y)@}_{\boldsymbol{\mathcal{B}}^{(2)}, \varphi(s)}\left(
\left\llbracket
\mathrm{CH}^{(2)}_{\boldsymbol{\mathcal{B}}^{(2)}, \varphi(s)}\left(
f^{(2)\flat}_{s}\left(
\mathrm{ech}^{(2,\mathcal{A}^{(1)})}_{\boldsymbol{\mathcal{A}}^{(2)}, s}\left(
\mathfrak{p}
\right)
\right)
\right)
\right\rrbracket_{\varphi(s)}
\right)
\tag{1}
\\
&=
\left\llbracket
\mathrm{ip}^{(2, Y)@}_{\boldsymbol{\mathcal{B}}^{(2)}, \varphi(s)}\left(
\mathrm{CH}^{(2)}_{\boldsymbol{\mathcal{B}}^{(2)}, \varphi(s)}\left(
f^{(2)\flat}_{s}\left(
\mathrm{ech}^{(2,\mathcal{A}^{(1)})}_{\boldsymbol{\mathcal{A}}^{(2)}, s}\left(
\mathfrak{p}
\right)
\right)
\right)
\right)
\right\rrbracket_{\varphi(s)}
\tag{2}
\\
&=
\left\llbracket
f^{(2)\flat}_{s}\left(
\mathrm{ech}^{(2,\mathcal{A}^{(1)})}_{\boldsymbol{\mathcal{A}}^{(2)}, s}\left(
\mathfrak{p}
\right)
\right)
\right\rrbracket_{\varphi(s)}
\tag{3}
\\
&=
\mathfrak{p}^{\llbracket\mathbf{Pth}_{\boldsymbol{\mathcal{B}}^{(2)}}^{\mathbf{f}^{(2)}}\rrbracket}
\tag{4}
\end{align*}
\end{flushleft}

The first equality unravels the definition of the constant $\mathfrak{p}$ in the partial $\Sigma^{\boldsymbol{\mathcal{A}}^{(2)}}$-algebra $\llbracket\mathbf{PT}_{\boldsymbol{\mathcal{B}}^{(2)}}^{\mathbf{f}^{(2)}}\rrbracket$, introduced in Proposition~\ref{PDQPTBDCatAlg};
the second equality unravels the definition of the mapping $\mathrm{ip}^{(\llbracket2\rrbracket, Y)@}_{\boldsymbol{\mathcal{B}}^{(2)}}$, introduced in Definition~\ref{DDPTQIp};
the third equality follows from Proposition~\ref{PDIpDCH};
finally, the last equality recovers the definition of the constant $\mathfrak{p}$ in the partial $\Sigma^{\boldsymbol{\mathcal{A}}^{(2)}}$-algebra $\llbracket\mathbf{Pth}_{\boldsymbol{\mathcal{B}}^{(2)}}^{\mathbf{f}^{(2)}}\rrbracket$, introduced in Proposition~\ref{PDQPthBDCatAlg}.

Hence, $\mathrm{ip}^{(\llbracket2\rrbracket, Y)@}_{\boldsymbol{\mathcal{B}}^{(2)}, \varphi}$ is compatible with the first-order rewrite rules.

{\sffamily The mapping $\mathrm{ip}^{(\llbracket2\rrbracket, Y)@}_{\boldsymbol{\mathcal{B}}^{(2)}, \varphi}$ is compatible with the $0$-source.}

Let $s$ be a sort in $S$ and let us consider the $0$-source operation symbol $\mathrm{sc}_{s}^{0}$ in $\Sigma^{\boldsymbol{\mathcal{A}}^{(2)}}_{s,s}$. Let $P$ be a second-order path term in $\mathrm{PT}_{\boldsymbol{\mathcal{B}}^{(2)}, \varphi(s)}$.

\allowdisplaybreaks
\begin{align*}
\mathrm{ip}^{(\llbracket2\rrbracket, Y)@}_{\boldsymbol{\mathcal{B}}^{(2)}, \varphi(s)} \left(
\mathrm{sc}_{s}^{0\llbracket\mathbf{PT}_{\boldsymbol{\mathcal{B}}^{(2)}}^{\mathbf{f}^{(2)}}\rrbracket}\left(
\left\llbracket
P
\right\rrbracket_{\varphi(s)}
\right)
\right)
&=
\mathrm{ip}^{(\llbracket2\rrbracket, Y)@}_{\boldsymbol{\mathcal{B}}^{(2)}, \varphi(s)} \left(
\mathrm{sc}_{\varphi(s)}^{0\llbracket\mathbf{PT}_{\boldsymbol{\mathcal{B}}^{(2)}}\rrbracket}\left(
\left\llbracket
P
\right\rrbracket_{\varphi(s)}
\right)
\right)
\tag{1}
\\
&=
\mathrm{sc}_{\varphi(s)}^{0\llbracket\mathbf{Pth}_{\boldsymbol{\mathcal{B}}^{(2)}}\rrbracket}\left(
\mathrm{ip}^{(\llbracket2\rrbracket, Y)@}_{\boldsymbol{\mathcal{B}}^{(2)}, \varphi(s)} \left(
\left\llbracket
P
\right\rrbracket_{\varphi(s)}
\right)
\right)
\tag{2}
\\
&=
\mathrm{sc}_{s}^{0\llbracket\mathbf{Pth}_{\boldsymbol{\mathcal{B}}^{(2)}}^{\mathbf{f}^{(2)}}\rrbracket}\left(
\mathrm{ip}^{(\llbracket2\rrbracket, Y)@}_{\boldsymbol{\mathcal{B}}^{(2)}, \varphi(s)} \left(
\left\llbracket
P
\right\rrbracket_{\varphi(s)}
\right)
\right).
\tag{3}
\end{align*}
The first equality unravels the interpretation of the operation symbol $\mathrm{sc}_{s}^{0}$ in the partial $\Sigma^{\boldsymbol{\mathcal{A}}^{(1)}}$-algebra $\llbracket\mathbf{PT}_{\boldsymbol{\mathcal{B}}^{(2)}}^{\mathbf{f}^{(2)}}\rrbracket$, introduced in Proposition~\ref{PDPTBDCatAlg};
the second equality follows from the fact that, according to Claim~\ref{CDIsoIpfc}, $\mathrm{ip}^{(\llbracket2\rrbracket, Y)@}_{\boldsymbol{\mathcal{B}}^{(2)}}$ is a $\Lambda^{\boldsymbol{\mathcal{B}}^{(2)}}$-homomorphism from $\llbracket\mathbf{PT}_{\boldsymbol{\mathcal{B}}^{(2)}}\rrbracket$ to $\llbracket\mathbf{Pth}_{\boldsymbol{\mathcal{B}}^{(2)}}\rrbracket$;
finally, the last equality recovers the interpretation of the operation symbol $\mathrm{sc}_{s}^{0}$ in the partial $\Sigma^{\boldsymbol{\mathcal{A}}^{(2)}}$-algebra $\llbracket\mathbf{Pth}_{\boldsymbol{\mathcal{B}}^{(2)}}^{\mathbf{f}^{(2)}}\rrbracket$, introduced in Proposition~\ref{PDQPthBCatAlg}.

Hence, $\mathrm{ip}^{(\llbracket2\rrbracket, Y)@}_{\boldsymbol{\mathcal{B}}^{(2)}, \varphi}$ is compatible with the $0$-source operation.

{\sffamily The mapping $\mathrm{ip}^{(\llbracket2\rrbracket, Y)@}_{\boldsymbol{\mathcal{B}}^{(2)}, \varphi}$ is compatible with the $0$-target.}

Let $s$ be a sort in $S$ and let us consider the $0$-target operation symbol $\mathrm{tg}_{s}^{0}$ in $\Sigma^{\boldsymbol{\mathcal{A}}^{(2)}}_{s,s}$. Let $P$ be a second-order path term in $\mathrm{PT}_{\boldsymbol{\mathcal{B}}^{(2)}, \varphi(s)}$, then the following equality holds
$$
\mathrm{ip}^{(\llbracket2\rrbracket, Y)@}_{\boldsymbol{\mathcal{B}}^{(2)}, \varphi(s)} \left(
\mathrm{tg}_{s}^{0\llbracket\mathbf{PT}_{\boldsymbol{\mathcal{B}}^{(2)}}^{\mathbf{f}^{(2)}}\rrbracket}\left(
\left\llbracket
P
\right\rrbracket_{\varphi(s)}
\right)
\right)
=
\mathrm{tg}_{s}^{0\llbracket\mathbf{Pth}_{\boldsymbol{\mathcal{B}}^{(2)}}^{\mathbf{f}^{(2)}}\rrbracket}\left(
\mathrm{ip}^{(\llbracket2\rrbracket, Y)@}_{\boldsymbol{\mathcal{B}}^{(2)}, \varphi(s)} \left(
\left\llbracket
P
\right\rrbracket_{\varphi(s)}
\right)
\right).
$$

The proof of this case is identical to that of the $0$-source.

Hence, $\mathrm{ip}^{(\llbracket2\rrbracket, Y)@}_{\boldsymbol{\mathcal{B}}^{(2)}, \varphi}$ is compatible with the $0$-target operation.

{\sffamily The mapping $\mathrm{ip}^{(\llbracket2\rrbracket, Y)@}_{\boldsymbol{\mathcal{B}}^{(2)}, \varphi}$ is compatible with the $0$-composition.}

Let $s$ be a sort in $S$ and let us consider the $0$-composition operation symbol $\circ_{s}^{0}$ in $\Sigma_{ss,s}^{\boldsymbol{\mathcal{A}}^{(2)}}$. Let $P$ and $Q$ be two second-order path terms in $\mathrm{PT}_{\boldsymbol{\mathcal{B}}^{(2)}, \varphi(s)}$ such that
$$
\mathrm{sc}_{\boldsymbol{\mathcal{B}}^{(2)}, \varphi(s)}^{(1,2)}\left(
\mathrm{ip}^{(2,Y)@}_{\boldsymbol{\mathcal{B}}^{(2)}, \varphi(s)}\left(
Q
\right)
\right)
=
\mathrm{tg}_{\boldsymbol{\mathcal{B}}^{(2)}, \varphi(s)}^{(1,2)}\left(
\mathrm{ip}^{(2,Y)@}_{\boldsymbol{\mathcal{B}}^{(2)}, \varphi(s)}\left(
P
\right)
\right).
$$
Then the following equality holds
\begin{multline*}
\mathrm{ip}^{(\llbracket2\rrbracket, Y)@}_{\boldsymbol{\mathcal{B}}^{(2)}, \varphi(s)} \left(
\left\llbracket
Q
\right\rrbracket_{\varphi(s)}
\circ_{s}^{0\llbracket\mathbf{PT}_{\boldsymbol{\mathcal{B}}^{(2)}}^{\mathbf{f}^{(2)}}\rrbracket}
\left\llbracket
P
\right\rrbracket_{\varphi(s)}
\right)
\\
=
\mathrm{ip}^{(\llbracket2\rrbracket, Y)@}_{\boldsymbol{\mathcal{B}}^{(2)}, \varphi(s)} \left(
\left\llbracket
Q
\right\rrbracket_{\varphi(s)}
\right)
\circ_{s}^{0\llbracket\mathbf{Pth}_{\boldsymbol{\mathcal{B}}^{(2)}}^{\mathbf{f}^{(2)}}\rrbracket}
\mathrm{ip}^{(\llbracket2\rrbracket, Y)@}_{\boldsymbol{\mathcal{B}}^{(2)}, \varphi(s)} \left(
\left\llbracket
P
\right\rrbracket_{\varphi(s)}
\right).
\end{multline*}

The proof of this case is identical to that of the $0$-source.

Hence, $\mathrm{ip}^{(\llbracket2\rrbracket, Y)@}_{\boldsymbol{\mathcal{B}}^{(2)}, \varphi}$ is compatible with the $0$-composition operation.

{\sffamily The mapping $\mathrm{ip}^{(\llbracket2\rrbracket, Y)@}_{\boldsymbol{\mathcal{B}}^{(2)}, \varphi}$ is compatible with the second-order rewrite rules.}

Let $s$ be a sort in $S$ and $\mathfrak{p}^{(2)}$ a rewrite rule in $\mathcal{A}_{s}^{(2)}$. Thus,
\begin{align*}
\mathrm{ip}^{(\llbracket2\rrbracket, Y)@}_{\boldsymbol{\mathcal{B}}^{(2)}, \varphi(s)}\left(
\mathfrak{p}^{(2)\llbracket\mathbf{PT}_{\boldsymbol{\mathcal{B}}^{(2)}}^{\mathbf{f}^{(2)}}\rrbracket}
\right)
&=
\mathrm{ip}^{(\llbracket2\rrbracket, Y)@}_{\boldsymbol{\mathcal{B}}^{(2)}, \varphi(s)}\left(
\left\llbracket
\mathrm{CH}^{(2)}_{\boldsymbol{\mathcal{B}}^{(2)}, \varphi(s)}\left(
f^{(2)\flat}_{s}\left(
\mathfrak{p}^{(2)}
\right)
\right)
\right\rrbracket_{\varphi(s)}
\right)
\tag{1}
\\
&=
\left\llbracket
\mathrm{ip}^{(2, Y)@}_{\boldsymbol{\mathcal{B}}^{(2)}, \varphi(s)}\left(
\mathrm{CH}^{(2)}_{\boldsymbol{\mathcal{B}}^{(2)}, \varphi(s)}\left(
f^{(2)\flat}_{s}\left(
\mathfrak{p}^{(2)}
\right)
\right)
\right)
\right\rrbracket_{\varphi(s)}
\tag{2}
\\
&=
\left\llbracket
f^{(2)\flat}_{s}\left(
\mathfrak{p}^{(2)}
\right)
\right\rrbracket_{\varphi(s)}
\tag{3}
\\
&=
\mathfrak{p}^{(2)\llbracket\mathbf{Pth}_{\boldsymbol{\mathcal{B}}^{(2)}}^{\mathbf{f}^{(2)}}\rrbracket}
\tag{4}
\end{align*}
The first equality unravels the definition of the constant $\mathfrak{p}^{(2)}$ in the partial $\Sigma^{\boldsymbol{\mathcal{A}}^{(2)}}$-algebra $\llbracket\mathbf{PT}_{\boldsymbol{\mathcal{B}}^{(2)}}^{\mathbf{f}^{(2)}}\rrbracket$, introduced in Proposition~\ref{PDQPTBDCatAlg};
the second equality unravels the definition of the mapping $\mathrm{ip}^{(\llbracket2\rrbracket, Y)@}_{\boldsymbol{\mathcal{B}}^{(2)}}$, introduced in Definition~\ref{DDPTQIp};
the third equality follows from Proposition~\ref{PDIpDCH};
finally, the last equality recovers the definition of the constant $\mathfrak{p}^{(2)}$ in the partial $\Sigma^{\boldsymbol{\mathcal{A}}^{(2)}}$-algebra $\llbracket\mathbf{Pth}_{\boldsymbol{\mathcal{B}}^{(2)}}^{\mathbf{f}^{(2)}}\rrbracket$, introduced in Proposition~\ref{PDQPthBDCatAlg}.

Hence, $\mathrm{ip}^{(\llbracket2\rrbracket, Y)@}_{\boldsymbol{\mathcal{B}}^{(2)}, \varphi}$ is compatible with the second-order rewrite rules.

{\sffamily The mapping $\mathrm{ip}^{(\llbracket2\rrbracket, Y)@}_{\boldsymbol{\mathcal{B}}^{(2)}, \varphi}$ is compatible with the $1$-source.}

Let $s$ be a sort in $S$ and let us consider the $1$-source operation symbol $\mathrm{sc}_{s}^{1}$ in $\Sigma^{\boldsymbol{\mathcal{A}}^{(2)}}_{s,s}$. Let $P$ be a second-order path term in $\mathrm{PT}_{\boldsymbol{\mathcal{B}}^{(2)}, \varphi(s)}$.

The following chain of equalities holds

\allowdisplaybreaks
\begin{align*}
\mathrm{ip}^{(\llbracket2\rrbracket, Y)@}_{\boldsymbol{\mathcal{B}}^{(2)}, \varphi(s)} \left(
\mathrm{sc}_{s}^{1\llbracket\mathbf{PT}_{\boldsymbol{\mathcal{B}}^{(2)}}^{\mathbf{f}^{(2)}}\rrbracket}\left(
\left\llbracket
P
\right\rrbracket_{\varphi(s)}
\right)
\right)
&=
\mathrm{ip}^{(\llbracket2\rrbracket, Y)@}_{\boldsymbol{\mathcal{B}}^{(2)}, \varphi(s)} \left(
\mathrm{sc}_{\varphi(s)}^{1\llbracket\mathbf{PT}_{\boldsymbol{\mathcal{B}}^{(2)}}\rrbracket}\left(
\left\llbracket
P
\right\rrbracket_{\varphi(s)}
\right)
\right)
\tag{1}
\\
&=
\mathrm{sc}_{\varphi(s)}^{1\llbracket\mathbf{Pth}_{\boldsymbol{\mathcal{B}}^{(2)}}\rrbracket}\left(
\mathrm{ip}^{(\llbracket2\rrbracket, Y)@}_{\boldsymbol{\mathcal{B}}^{(2)}, \varphi(s)} \left(
\left\llbracket
P
\right\rrbracket_{\varphi(s)}
\right)
\right)
\tag{2}
\\
&=
\mathrm{sc}_{s}^{1\llbracket\mathbf{Pth}_{\boldsymbol{\mathcal{B}}^{(2)}}^{\mathbf{f}^{(2)}}\rrbracket}\left(
\mathrm{ip}^{(\llbracket2\rrbracket, Y)@}_{\boldsymbol{\mathcal{B}}^{(2)}, \varphi(s)} \left(
\left\llbracket
P
\right\rrbracket_{\varphi(s)}
\right)
\right).
\tag{3}
\end{align*}
The first equality unravels the interpretation of the operation symbol $\mathrm{sc}_{s}^{1}$ in the partial $\Sigma^{\boldsymbol{\mathcal{A}}^{(1)}}$-algebra $\llbracket\mathbf{PT}_{\boldsymbol{\mathcal{B}}^{(2)}}^{\mathbf{f}^{(2)}}\rrbracket$, introduced in Proposition~\ref{PDPTBDCatAlg};
the second equality follows from the fact that, according to Claim~\ref{CDIsoIpfc}, $\mathrm{ip}^{(\llbracket2\rrbracket, X)@}_{\boldsymbol{\mathcal{B}}^{(2)}}$ is a $\Lambda^{\boldsymbol{\mathcal{B}}^{(2)}}$-homomorphism from $\llbracket\mathbf{PT}_{\boldsymbol{\mathcal{B}}^{(2)}}\rrbracket$ to $\llbracket\mathbf{Pth}_{\boldsymbol{\mathcal{B}}^{(2)}}\rrbracket$;
finally, the last equality recovers the interpretation of the operation symbol $\mathrm{sc}_{s}^{1}$ in the partial $\Sigma^{\boldsymbol{\mathcal{A}}^{(2)}}$-algebra $\llbracket\mathbf{Pth}_{\boldsymbol{\mathcal{B}}^{(2)}}^{\mathbf{f}^{(2)}}\rrbracket$, introduced in Proposition~\ref{PDQPthBCatAlg}.

Hence, $\mathrm{ip}^{(\llbracket2\rrbracket, Y)@}_{\boldsymbol{\mathcal{B}}^{(2)}, \varphi}$ is compatible with the $1$-source operation.

{\sffamily The mapping $\mathrm{ip}^{(\llbracket2\rrbracket, Y)@}_{\boldsymbol{\mathcal{B}}^{(2)}, \varphi}$ is compatible with the $1$-target.}

Let $s$ be a sort in $S$ and let us consider the $1$-target operation symbol $\mathrm{tg}_{s}^{1}$ in $\Sigma^{\boldsymbol{\mathcal{A}}^{(2)}}_{s,s}$. Let $P$ be a second-order path term in $\mathrm{PT}_{\boldsymbol{\mathcal{B}}^{(2)}, \varphi(s)}$, then the following equality holds
$$
\mathrm{ip}^{(\llbracket2\rrbracket, Y)@}_{\boldsymbol{\mathcal{B}}^{(2)}, \varphi(s)} \left(
\mathrm{tg}_{s}^{1\llbracket\mathbf{PT}_{\boldsymbol{\mathcal{B}}^{(2)}}^{\mathbf{f}^{(2)}}\rrbracket}\left(
\left\llbracket
P
\right\rrbracket_{\varphi(s)}
\right)
\right)
=
\mathrm{tg}_{s}^{1\llbracket\mathbf{Pth}_{\boldsymbol{\mathcal{B}}^{(2)}}^{\mathbf{f}^{(2)}}\rrbracket}\left(
\mathrm{ip}^{(\llbracket2\rrbracket, Y)@}_{\boldsymbol{\mathcal{B}}^{(2)}, \varphi(s)} \left(
\left\llbracket
P
\right\rrbracket_{\varphi(s)}
\right)
\right).
$$

The proof of this case is identical to that of the $1$-source.

Hence, $\mathrm{ip}^{(\llbracket2\rrbracket, Y)@}_{\boldsymbol{\mathcal{B}}^{(2)}, \varphi}$ is compatible with the $1$-target operation.

{\sffamily The mapping $\mathrm{ip}^{(\llbracket2\rrbracket, Y)@}_{\boldsymbol{\mathcal{B}}^{(2)}, \varphi}$ is compatible with the $1$-composition.}

Let $s$ be a sort in $S$ and let us consider the $1$-composition operation symbol $\circ_{s}^{1}$ in $\Sigma_{ss,s}^{\boldsymbol{\mathcal{A}}^{(2)}}$. Let $P$ and $Q$ be two second-order path terms in $\mathrm{PT}_{\boldsymbol{\mathcal{B}}^{(2)}, \varphi(s)}$ such that
$$
\mathrm{sc}_{\boldsymbol{\mathcal{B}}^{(2)}, \varphi(s)}^{([1],2)}\left(
\mathrm{ip}^{(2,Y)@}_{\boldsymbol{\mathcal{B}}^{(2)}, \varphi(s)}\left(
Q
\right)
\right)
=
\mathrm{tg}_{\boldsymbol{\mathcal{B}}^{(2)}, \varphi(s)}^{([1],2)}\left(
\mathrm{ip}^{(2,Y)@}_{\boldsymbol{\mathcal{B}}^{(2)}, \varphi(s)}\left(
P
\right)
\right).
$$
Then the following equality holds
\begin{multline*}
\mathrm{ip}^{(\llbracket2\rrbracket, Y)@}_{\boldsymbol{\mathcal{B}}^{(2)}, \varphi(s)} \left(
\left\llbracket
Q
\right\rrbracket_{\varphi(s)}
\circ_{s}^{1\llbracket\mathbf{PT}_{\boldsymbol{\mathcal{B}}^{(2)}}^{\mathbf{f}^{(2)}}\rrbracket}
\left\llbracket
P
\right\rrbracket_{\varphi(s)}
\right)
\\
=
\mathrm{ip}^{(\llbracket2\rrbracket, Y)@}_{\boldsymbol{\mathcal{B}}^{(2)}, \varphi(s)} \left(
\left\llbracket
Q
\right\rrbracket_{\varphi(s)}
\right)
\circ_{s}^{1\llbracket\mathbf{Pth}_{\boldsymbol{\mathcal{B}}^{(2)}}^{\mathbf{f}^{(2)}}\rrbracket}
\mathrm{ip}^{(\llbracket2\rrbracket, Y)@}_{\boldsymbol{\mathcal{B}}^{(2)}, \varphi(s)} \left(
\left\llbracket
P
\right\rrbracket_{\varphi(s)}
\right).
\end{multline*}

The proof of this case is identical to that of the $1$-source.

Hence, $\mathrm{ip}^{(\llbracket2\rrbracket, Y)@}_{\boldsymbol{\mathcal{B}}^{(2)}, \varphi}$ is compatible with the $1$-composition operation.

This completes the proof.
\end{proof}

\begin{theorem}\label{TDIsoB}
The partial $\Sigma^{\boldsymbol{\mathcal{A}}^{(2)}}$-algebras $\llbracket\mathbf{Pth}_{\boldsymbol{\mathcal{B}}^{(2)}}^{\mathbf{f}^{(2)}}\rrbracket$ and $\llbracket\mathbf{PT}_{\boldsymbol{\mathcal{B}}^{(2)}}^{\mathbf{f}^{(2)}}\rrbracket$ are isomorphic.
\end{theorem}

\begin{proof}
Let $s$ be a sort in $S$ and $\mathfrak{P}^{(2)}$ be a second-order path in $\mathrm{Pth}_{\boldsymbol{\mathcal{B}}^{(2)}, \varphi(s)}$.

The following chain of equalities holds
\allowdisplaybreaks
\begin{flushleft}
$
\mathrm{ip}^{(\llbracket2\rrbracket, Y)@}_{\boldsymbol{\mathcal{B}}^{(2)}, \varphi(s)} \left(
\mathrm{CH}^{\llbracket2\rrbracket}_{\boldsymbol{\mathcal{B}}^{(2)}, \varphi(s)} \left(
\left\llbracket
\mathfrak{P}^{(2)}
\right\rrbracket_{\varphi(s)}
\right)
\right)
$
\begin{align*}
&=
\left\llbracket
\mathrm{ip}^{(2, Y)@}_{\boldsymbol{\mathcal{B}}^{(2)}, \varphi(s)}\left(
\mathrm{CH}^{(2)}_{\boldsymbol{\mathcal{B}}^{(2)}, \varphi(s)}\left(
\mathfrak{P}^{(2)}
\right)
\right)
\right\rrbracket_{\varphi(s)}
\tag{1}
\\
&=
\left\llbracket
\mathfrak{P}^{(2)}
\right\rrbracket_{\varphi(s)}.
\tag{2}
\end{align*}
\end{flushleft}

The first equality unravels the definition of the mappings $\mathrm{ip}^{(\llbracket2\rrbracket, Y)@}_{\boldsymbol{\mathcal{B}}^{(2)}, \varphi}$ and $\mathrm{CH}^{\llbracket2\rrbracket}_{\boldsymbol{\mathcal{B}}^{(2)}, \varphi}$ according to, respectively, Definitions~\ref{DDPTQIp} and \ref{DDPTQDCH};
finally, the second equality follows from Proposition~\ref{PDIpDCH}.

On the other hand, let $s$ be a sort in $S$ and $P$ a second-order path term in $\mathrm{PT}_{\boldsymbol{\mathcal{B}}^{(2)}, \varphi(s)}$.

The following chain of equalities holds
\allowdisplaybreaks
\begin{align*}
\mathrm{CH}^{\llbracket2\rrbracket}_{\boldsymbol{\mathcal{B}}^{(2)}, \varphi(s)} \left(
\mathrm{ip}^{(\llbracket2\rrbracket, Y)@}_{\boldsymbol{\mathcal{B}}^{(2)}, \varphi(s)} \left(
\left\llbracket
P
\right\rrbracket_{\varphi(s)}
\right)
\right)
&=
\left\llbracket
\mathrm{CH}^{(2)}_{\boldsymbol{\mathcal{B}}^{(2)}, \varphi(s)}\left(
\mathrm{ip}^{(2, Y)@}_{\boldsymbol{\mathcal{B}}^{(2)}, \varphi(s)}\left(
P
\right)
\right)
\right\rrbracket_{\varphi(s)}
\tag{1}
\\
&=
\left\llbracket
P
\right\rrbracket_{\varphi(s)}.
\tag{2}
\end{align*}
The first equality unravels the definition of the mappings $\mathrm{ip}^{(\llbracket2\rrbracket, Y)@}_{\boldsymbol{\mathcal{B}}^{(2)}, \varphi}$ and $\mathrm{CH}^{\llbracket2\rrbracket}_{\boldsymbol{\mathcal{B}}^{(2)}, \varphi}$ according to, respectively, Definitions~\ref{DDPTQIp} and \ref{DDPTQDCH};
finally, the second equality follows from Proposition~\ref{PDIpDCH}.
\end{proof}

\section{A structure of partial $\Sigma^{\boldsymbol{\mathcal{A}}^{(2)}}$-algebra on $\T_{\boldsymbol{\mathcal{E}}^{\boldsymbol{\mathcal{B}}^{(2)}}}(\mathbf{Pth}_{\boldsymbol{\mathcal{B}}^{(2)}})$}

We show that $\T_{\boldsymbol{\mathcal{E}}^{\boldsymbol{\mathcal{B}}^{(2)}}}(\mathbf{Pth}_{\boldsymbol{\mathcal{B}}^{(2)}})_{\varphi}$ is equipped, in a natural way with a structure of partial $\Sigma^{\boldsymbol{\mathcal{A}}^{(2)}}$-algebra.

\begin{proposition}
\label{PDFreeBDCatAlg}
Let $\mathbf{f}^{(2)}=(\varphi, c, (f^{(i)})_{i\in 3})$ be a second-order morphism from $\boldsymbol{\mathcal{A}}^{(2)}$ to $\boldsymbol{\mathcal{B}}^{(2)}$. Then the $S$-sorted set $\T_{\boldsymbol{\mathcal{E}}^{\boldsymbol{\mathcal{B}}^{(2)}}}(\mathbf{Pth}_{\boldsymbol{\mathcal{B}}^{(2)}})_{\varphi}$ is equipped, in a natural way, with a structure of partial $\Sigma^{\boldsymbol{\mathcal{A}}^{(2)}}$-algebra.
\end{proposition}

\begin{proof}
Let us denote by $\mathbf{T}_{\boldsymbol{\mathcal{E}}^{\boldsymbol{\mathcal{B}}^{(2)}}}^{\mathbf{f}^{(2)}}(\mathbf{Pth}_{\boldsymbol{\mathcal{B}}^{(2)}})$ the $\Sigma^{\boldsymbol{\mathcal{A}}^{(2)}}$-algebra defined as follows:

\textsf{(1)}
The underlying $S$-sorted set of $\mathbf{T}_{\boldsymbol{\mathcal{E}}^{\boldsymbol{\mathcal{B}}^{(2)}}}^{\mathbf{f}^{(2)}}(\mathbf{Pth}_{\boldsymbol{\mathcal{B}}^{(2)}})$ is $\T_{\boldsymbol{\mathcal{E}}^{\boldsymbol{\mathcal{B}}^{(2)}}}(\mathbf{Pth}_{\boldsymbol{\mathcal{B}}^{(2)}})_{\varphi}=(\T_{\boldsymbol{\mathcal{E}}^{\boldsymbol{\mathcal{B}}^{(2)}}}(\mathbf{Pth}_{\boldsymbol{\mathcal{B}}^{(2)}})_{\varphi(s)})_{s\in S}$.

\textsf{(2)}
For every $(\mathbf{s}, s)\in S^{\ast}\times S$ and every operation symbol $\sigma\in\Sigma_{\mathbf{s},s}$, the operation $\sigma^{\mathbf{T}_{\boldsymbol{\mathcal{E}}^{\boldsymbol{\mathcal{B}}^{(2)}}}^{\mathbf{f}^{(2)}}(\mathbf{Pth}_{\boldsymbol{\mathcal{B}}^{(2)}})}$ is given by the interpretation of $\sigma$ in the $\Sigma$-algebra $\mathbf{c}_{\mathfrak{d}}^{\ast}(\mathbf{T}_{\boldsymbol{\mathcal{E}}^{\boldsymbol{\mathcal{B}}^{(2)}}}^{(0,2)}(\mathbf{Pth}_{\boldsymbol{\mathcal{B}}^{(2)}}))$. That is, its interpretation is given by the derived operation in $\mathbf{T}_{\boldsymbol{\mathcal{E}}^{\boldsymbol{\mathcal{B}}^{(2)}}}^{(0,2)}(\mathbf{Pth}_{\boldsymbol{\mathcal{B}}^{(2)}})$, the $\Lambda$-reduct of the $\Lambda^{\boldsymbol{\mathcal{B}}^{(2)}}$-algebra  $\mathbf{T}_{\boldsymbol{\mathcal{E}}^{\boldsymbol{\mathcal{B}}^{(2)}}}(\mathbf{Pth}_{\boldsymbol{\mathcal{B}}^{(2)}})$, introduced in Definition~\ref{DDVarAbv}.

\textsf{(3)}
For every $s\in S$ and every $\mathfrak{p}\in\mathcal{A}_{s}^{(1)}$, the constant $\mathfrak{p}^{\mathbf{T}_{\boldsymbol{\mathcal{E}}^{\boldsymbol{\mathcal{B}}^{(2)}}}^{\mathbf{f}^{(2)}}(\mathbf{Pth}_{\boldsymbol{\mathcal{B}}^{(2)}})}$ is given by
$$
\mathfrak{p}^{\mathbf{T}_{\boldsymbol{\mathcal{E}}^{\boldsymbol{\mathcal{B}}^{(2)}}}^{\mathbf{f}^{(2)}}(\mathbf{Pth}_{\boldsymbol{\mathcal{B}}^{(2)}})}
=
\mathrm{pr}^{\equiv^{\llbracket2\rrbracket}}_{\boldsymbol{\mathcal{B}}^{(2)}, \varphi(s)} \circ \mathrm{ip}^{(2,Y)@}_{\boldsymbol{\mathcal{B}}^{(2)}, \varphi(s)} \circ \mathrm{CH}^{(2)}_{\boldsymbol{\mathcal{B}}^{(2)}, \varphi(s)} \left(
f^{(2)\flat}_{s}\left(
\mathfrak{p}^{\mathbf{Pth}_{\boldsymbol{\mathcal{A}}^{(2)}}}
\right)
\right).
$$

\textsf{(4)}
For every $s\in S$, the interpretations of the operations $\mathrm{sc}_{s}^{0}$ and $\mathrm{tg}_{s}^{0}$ are given by
\begin{align*}
\mathrm{sc}_{s}^{0\mathbf{T}_{\boldsymbol{\mathcal{E}}^{\boldsymbol{\mathcal{B}}^{(2)}}}^{\mathbf{f}^{(2)}}(\mathbf{Pth}_{\boldsymbol{\mathcal{B}}^{(2)}})}
&=
\mathrm{sc}_{\varphi(s)}^{0\mathbf{T}_{\boldsymbol{\mathcal{E}}^{\boldsymbol{\mathcal{B}}^{(2)}}}(\mathbf{Pth}_{\boldsymbol{\mathcal{B}}^{(2)}})}
&&\mbox{and}&
\mathrm{tg}_{s}^{0\mathbf{T}_{\boldsymbol{\mathcal{E}}^{\boldsymbol{\mathcal{B}}^{(2)}}}^{\mathbf{f}^{(2)}}(\mathbf{Pth}_{\boldsymbol{\mathcal{B}}^{(2)}})}
&=
\mathrm{tg}_{\varphi(s)}^{0\mathbf{T}_{\boldsymbol{\mathcal{E}}^{\boldsymbol{\mathcal{B}}^{(2)}}}(\mathbf{Pth}_{\boldsymbol{\mathcal{B}}^{(2)}})}.
\end{align*}
That is, their interpretations are given by the interpretations of $\mathrm{sc}_{\varphi(s)}^{0}$ and $\mathrm{tg}_{\varphi(s)}^{0}$ in $\mathbf{T}_{\boldsymbol{\mathcal{E}}^{\boldsymbol{\mathcal{B}}^{(2)}}}(\mathbf{Pth}_{\boldsymbol{\mathcal{B}}^{(2)}})$ that, we recall, was introduced in Definition~\ref{DDVarAbv}.

\textsf{(5)}
Similarly,  for every $s \in S$, the interpretations of the operation $\circ_{s}^{0}$ is given by
$$
\circ_{s}^{0\mathbf{T}_{\boldsymbol{\mathcal{E}}^{\boldsymbol{\mathcal{B}}^{(2)}}}^{\mathbf{f}^{(2)}}(\mathbf{Pth}_{\boldsymbol{\mathcal{B}}^{(2)}})}
=
\circ_{\varphi(s)}^{0\mathbf{T}_{\boldsymbol{\mathcal{E}}^{\boldsymbol{\mathcal{B}}^{(2)}}}(\mathbf{Pth}_{\boldsymbol{\mathcal{B}}^{(2)}})}.
$$
That is, its interpretation is given by the interpretation of $\circ_{\varphi(s)}^{0}$ in $\mathbf{T}_{\boldsymbol{\mathcal{E}}^{\boldsymbol{\mathcal{B}}^{(2)}}}(\mathbf{Pth}_{\boldsymbol{\mathcal{B}}^{(2)}})$ that, we recall, was introduced in Definition~\ref{DDVarAbv}.

\textsf{(6)}
For every $s\in S$ and every $\mathfrak{p}\in\mathcal{A}_{s}^{(2)}$, the constant $\mathfrak{p}^{\mathbf{Pth}_{\boldsymbol{\mathcal{B}}^{(2)}}^{\mathbf{f}^{(2)}}}$ is given by
$$
\mathfrak{p}^{(2)\mathbf{PT}_{\boldsymbol{\mathcal{B}}^{(2)}}^{\mathbf{f}^{(2)}}}
=
\mathrm{pr}^{\equiv^{\llbracket2\rrbracket}}_{\boldsymbol{\mathcal{B}}^{(2)}, \varphi(s)} \circ \mathrm{ip}^{(2,Y)@}_{\boldsymbol{\mathcal{B}}^{(2)}, \varphi(s)} \circ \mathrm{CH}^{(2)}_{\boldsymbol{\mathcal{B}}^{(2)}, \varphi(s)} \left(
f^{(2)}_{s}\left(
\mathfrak{p}^{(2)}
\right)
\right).
$$

\textsf{(7)}
For every $s\in S$, the interpretations of the operations $\mathrm{sc}_{s}^{1}$ and $\mathrm{tg}_{s}^{1}$ are given by
\begin{align*}
\mathrm{sc}_{s}^{1\mathbf{T}_{\boldsymbol{\mathcal{E}}^{\boldsymbol{\mathcal{B}}^{(2)}}}^{\mathbf{f}^{(2)}}(\mathbf{Pth}_{\boldsymbol{\mathcal{B}}^{(2)}})}
&=
\mathrm{sc}_{\varphi(s)}^{1\mathbf{T}_{\boldsymbol{\mathcal{E}}^{\boldsymbol{\mathcal{B}}^{(2)}}}(\mathbf{Pth}_{\boldsymbol{\mathcal{B}}^{(2)}})}
&&\mbox{and}&
\mathrm{tg}_{s}^{1\mathbf{T}_{\boldsymbol{\mathcal{E}}^{\boldsymbol{\mathcal{B}}^{(2)}}}^{\mathbf{f}^{(2)}}(\mathbf{Pth}_{\boldsymbol{\mathcal{B}}^{(2)}})}
&=
\mathrm{tg}_{\varphi(s)}^{1\mathbf{T}_{\boldsymbol{\mathcal{E}}^{\boldsymbol{\mathcal{B}}^{(2)}}}(\mathbf{Pth}_{\boldsymbol{\mathcal{B}}^{(2)}})}.
\end{align*}
That is, their interpretations are given by the interpretations of $\mathrm{sc}_{\varphi(s)}^{1}$ and $\mathrm{tg}_{\varphi(s)}^{1}$ in $\mathbf{T}_{\boldsymbol{\mathcal{E}}^{\boldsymbol{\mathcal{B}}^{(2)}}}(\mathbf{Pth}_{\boldsymbol{\mathcal{B}}^{(2)}})$ that, we recall, was introduced in Definition~\ref{DDVarAbv}.

\textsf{(8)}
Similarly,  for every $s \in S$, the interpretations of the operation $\circ_{s}^{1}$ is given by
$$
\circ_{s}^{1\mathbf{T}_{\boldsymbol{\mathcal{E}}^{\boldsymbol{\mathcal{B}}^{(2)}}}^{\mathbf{f}^{(2)}}(\mathbf{Pth}_{\boldsymbol{\mathcal{B}}^{(2)}})}
=
\circ_{\varphi(s)}^{1\mathbf{T}_{\boldsymbol{\mathcal{E}}^{\boldsymbol{\mathcal{B}}^{(2)}}}(\mathbf{Pth}_{\boldsymbol{\mathcal{B}}^{(2)}})}.
$$
That is, its interpretation is given by the interpretation of $\circ_{\varphi(s)}^{1}$ in $\mathbf{T}_{\boldsymbol{\mathcal{E}}^{\boldsymbol{\mathcal{B}}^{(2)}}}(\mathbf{Pth}_{\boldsymbol{\mathcal{B}}^{(2)}})$ that, we recall, was introduced in Definition~\ref{DDVarAbv}.

This completes the definition of the partial $\Sigma^{\boldsymbol{\mathcal{A}}^{(2)}}$-algebra $\mathbf{T}_{\boldsymbol{\mathcal{E}}^{\boldsymbol{\mathcal{B}}^{(2)}}}^{\mathbf{f}^{(2)}}(\mathbf{Pth}_{\boldsymbol{\mathcal{B}}^{(2)}})$.
\end{proof}

Finally, we show that the partial $\Sigma^{\boldsymbol{\mathcal{A}}^{(2)}}$-algebras $\llbracket\mathbf{Pth}_{\boldsymbol{\mathcal{B}}^{(2)}}^{\mathbf{f}^{(2)}}\rrbracket_{\varphi}$ and $\mathbf{T}_{\boldsymbol{\mathcal{E}}^{\boldsymbol{\mathcal{B}}^{(2)}}}^{\mathbf{f}^{(2)}}(\mathbf{Pth}_{\boldsymbol{\mathcal{B}}^{(2)}})$ are isomorphic. In virtue of Theorem~\ref{TDIsoB} the same relation will apply to the partial $\Sigma^{\boldsymbol{\mathcal{A}}^{(2)}}$-algebra $\llbracket \mathbf{PT}^{\mathbf{f}^{(2)}}_{\boldsymbol{\mathcal{B}}^{(2)}} \rrbracket$.

\begin{proposition}
\label{PDQCHBDCatHom}
Let $\mathbf{f}^{(2)}=(\varphi, c, (f^{(i)})_{i\in 3})$ be a second-order morphism from $\boldsymbol{\mathcal{A}}^{(2)}$ to $\boldsymbol{\mathcal{B}}^{(2)}$. Then the mapping 
$$
(\mathrm{pr}^{\equiv^{\llbracket2\rrbracket}}_{\boldsymbol{\mathcal{B}}^{(2)}} \circ \mathrm{ip}^{(2,Y)@}_{\boldsymbol{\mathcal{B}}^{(2)}} \circ \mathrm{CH}^{(2)}_{\boldsymbol{\mathcal{B}}^{(2)}})^{\natural}_{\varphi}
\colon 
\llbracket\mathrm{Pth}_{\boldsymbol{\mathcal{B}}^{(2)}}\rrbracket_{\varphi}
\mor
\T_{\boldsymbol{\mathcal{E}}^{\boldsymbol{\mathcal{B}}^{(2)}}}(\mathbf{Pth}_{\boldsymbol{\mathcal{B}}^{(2)}})_{\varphi}
$$
is a $\Sigma^{\boldsymbol{\mathcal{A}}^{(2)}}$-homomorphism from $\llbracket\mathbf{Pth}_{\boldsymbol{\mathcal{B}}^{(2)}}^{\mathbf{f}^{(2)}}\rrbracket$ to $\mathbf{T}_{\boldsymbol{\mathcal{E}}^{\boldsymbol{\mathcal{B}}^{(2)}}}^{\mathbf{f}^{(2)}}(\mathbf{Pth}_{\boldsymbol{\mathcal{B}}^{(2)}})$.
\end{proposition}

\begin{proof}
We prove that $(\mathrm{pr}^{\equiv^{\llbracket2\rrbracket}}_{\boldsymbol{\mathcal{B}}^{(2)}} \circ \mathrm{ip}^{(2,Y)@}_{\boldsymbol{\mathcal{B}}^{(2)}} \circ \mathrm{CH}^{(2)}_{\boldsymbol{\mathcal{B}}^{(2)}})^{\natural}_{\varphi}$ is compatible with every operation symbol in $\Sigma^{\boldsymbol{\mathcal{A}}^{(2)}}$.

{\sffamily The mapping $(\mathrm{pr}^{\equiv^{\llbracket2\rrbracket}}_{\boldsymbol{\mathcal{B}}^{(2)}} \circ \mathrm{ip}^{(2,Y)@}_{\boldsymbol{\mathcal{B}}^{(2)}} \circ \mathrm{CH}^{(2)}_{\boldsymbol{\mathcal{B}}^{(2)}})^{\natural}_{\varphi}$ is a $\Sigma$-homomorphism.}

Note that $(\mathrm{pr}^{\equiv^{\llbracket2\rrbracket}}_{\boldsymbol{\mathcal{B}}^{(2)}} \circ \mathrm{ip}^{(2,Y)@}_{\boldsymbol{\mathcal{B}}^{(2)}} \circ \mathrm{CH}^{(2)}_{\boldsymbol{\mathcal{B}}^{(2)}})^{\natural}_{\varphi} = \mathbf{c}_{\mathfrak{d}}^{\ast} ((\mathrm{pr}^{\equiv^{\llbracket2\rrbracket}}_{\boldsymbol{\mathcal{B}}^{(2)}} \circ \mathrm{ip}^{(2,Y)@}_{\boldsymbol{\mathcal{B}}^{(2)}} \circ \mathrm{CH}^{(2)}_{\boldsymbol{\mathcal{B}}^{(2)}})^{\natural})$. By Remark~\ref{RDQPUniv}, the mapping $(\mathrm{pr}^{\equiv^{\llbracket2\rrbracket}}_{\boldsymbol{\mathcal{B}}^{(2)}} \circ \mathrm{ip}^{(2,Y)@}_{\boldsymbol{\mathcal{B}}^{(2)}} \circ \mathrm{CH}^{(2)}_{\boldsymbol{\mathcal{B}}^{(2)}})^{\natural}$ is a $\Lambda^{\boldsymbol{\mathcal{B}}^{(2)}}$-homomorphism, thus in particular a $\Lambda$-homomorphism. Therefore, it follows from Proposition~\ref{PFunSig} that the mapping $(\mathrm{pr}^{\equiv^{\llbracket2\rrbracket}}_{\boldsymbol{\mathcal{B}}^{(2)}} \circ \mathrm{ip}^{(2,Y)@}_{\boldsymbol{\mathcal{B}}^{(2)}} \circ \mathrm{CH}^{(2)}_{\boldsymbol{\mathcal{B}}^{(2)}})^{\natural}_{\varphi}$ is a $\Sigma$-homomorphism.


{\sffamily The mapping $(\mathrm{pr}^{\equiv^{\llbracket2\rrbracket}}_{\boldsymbol{\mathcal{B}}^{(2)}} \circ \mathrm{ip}^{(2,Y)@}_{\boldsymbol{\mathcal{B}}^{(2)}} \circ \mathrm{CH}^{(2)}_{\boldsymbol{\mathcal{B}}^{(2)}})^{\natural}_{\varphi}$ is compatible with the first-order rewrite rules.}

Let $s$ be a sort in $S$ and $\mathfrak{p}$ a rewrite rule in $\mathcal{A}_{s}^{(1)}$. Thus, the following chain of equalities holds
\begin{flushleft}
$
(\mathrm{pr}^{\equiv^{\llbracket2\rrbracket}}_{\boldsymbol{\mathcal{B}}^{(2)}} \circ \mathrm{ip}^{(2,Y)@}_{\boldsymbol{\mathcal{B}}^{(2)}} \circ \mathrm{CH}^{(2)}_{\boldsymbol{\mathcal{B}}^{(2)}})^{\natural}_{\varphi(s)} \left(
\mathfrak{p}^{\llbracket\mathbf{Pth}_{\boldsymbol{\mathcal{B}}^{(2)}}^{\mathbf{f}^{(2)}}\rrbracket}
\right)
$
\begin{align*}
&=
(\mathrm{pr}^{\equiv^{\llbracket2\rrbracket}}_{\boldsymbol{\mathcal{B}}^{(2)}} \circ \mathrm{ip}^{(2,Y)@}_{\boldsymbol{\mathcal{B}}^{(2)}} \circ \mathrm{CH}^{(2)}_{\boldsymbol{\mathcal{B}}^{(2)}})^{\natural}_{\varphi(s)} \left(
\left\llbracket
f^{(2)\flat}\left(
\mathfrak{p}^{\mathbf{Pth}_{\boldsymbol{\mathcal{A}}^{(2)}}}
\right)
\right\rrbracket_{\varphi(s)}
\right)
\tag{1}
\\
&=
(\mathrm{pr}^{\equiv^{\llbracket2\rrbracket}}_{\boldsymbol{\mathcal{B}}^{(2)}} \circ \mathrm{ip}^{(2,Y)@}_{\boldsymbol{\mathcal{B}}^{(2)}} \circ \mathrm{CH}^{(2)}_{\boldsymbol{\mathcal{B}}^{(2)}})_{\varphi(s)} \left(
f^{(2)\flat}\left(
\mathfrak{p}^{\mathbf{Pth}_{\boldsymbol{\mathcal{A}}^{(2)}}}
\right)
\right)
\tag{2}
\\
&=
\mathrm{pr}^{\equiv^{\llbracket2\rrbracket}}_{\boldsymbol{\mathcal{B}}^{(2)}, \varphi(s)} 
\circ
\mathrm{ip}^{(2,Y)@}_{\boldsymbol{\mathcal{B}}^{(2)}, \varphi(s)}
\circ
\mathrm{CH}^{(2)}_{\boldsymbol{\mathcal{B}}^{(2)}, \varphi(s)} \left(
f^{(2)\flat}\left(
\mathfrak{p}^{\mathbf{Pth}_{\boldsymbol{\mathcal{A}}^{(2)}}}
\right)
\right)
\tag{3}
\\
&=
\mathfrak{p}^{\mathbf{T}_{\boldsymbol{\mathcal{E}}^{\boldsymbol{\mathcal{B}}^{(2)}}}^{\mathbf{f}^{(2)}}(\mathbf{Pth}_{\boldsymbol{\mathcal{B}}^{(2)}})}.
\tag{4}
\end{align*}
\end{flushleft}

The first equality unravels the interpretation of the constant $\mathfrak{p}$ in the partial $\Sigma^{\boldsymbol{\mathcal{A}}^{(2)}}$-algebra $\llbracket\mathbf{Pth}^{\mathbf{f}^{(2)}}_{\boldsymbol{\mathcal{B}}^{(2)}}\rrbracket$, introduced in Proposition~\ref{PDQPthBDCatAlg};
the second equality follows from Remark~\ref{RDQPUniv};
the third equality unravels the definition of the $\varphi(s)$ component of a composition of many-sorted mappings;
finally, the last equality recovers the interpretation of the constant $\mathfrak{p}$ in the partial $\Sigma^{\boldsymbol{\mathcal{A}}^{(2)}}$-algebra $\mathbf{T}_{\boldsymbol{\mathcal{E}}^{\boldsymbol{\mathcal{B}}^{(2)}}}^{\mathbf{f}^{(2)}}(\mathbf{Pth}_{\boldsymbol{\mathcal{B}}^{(2)}})$, introduced in Proposition~\ref{PDFreeBDCatAlg}.

Hence, $(\mathrm{pr}^{\equiv^{\llbracket2\rrbracket}}_{\boldsymbol{\mathcal{B}}^{(2)}} \circ \mathrm{ip}^{(2,Y)@}_{\boldsymbol{\mathcal{B}}^{(2)}} \circ \mathrm{CH}^{(2)}_{\boldsymbol{\mathcal{B}}^{(2)}})^{\natural}_{\varphi}$ is compatible with the first-order rewrite rules.

{\sffamily The mapping $(\mathrm{pr}^{\equiv^{\llbracket2\rrbracket}}_{\boldsymbol{\mathcal{B}}^{(2)}} \circ \mathrm{ip}^{(2,Y)@}_{\boldsymbol{\mathcal{B}}^{(2)}} \circ \mathrm{CH}^{(2)}_{\boldsymbol{\mathcal{B}}^{(2)}})^{\natural}_{\varphi}$ is compatible with the $0$-source.}

Let $s$ be a sort in $S$ and let us consider the $0$-source operation symbol $\mathrm{sc}_{s}^{0}$ in $\Sigma^{\boldsymbol{\mathcal{A}}^{(1)}}_{s,s}$. Let $\mathfrak{P}^{(2)}$ be a second-order path in $\mathrm{Pth}_{\boldsymbol{\mathcal{B}}^{(2)}, \varphi(s)}$.

The following chain of equalities holds
\begin{flushleft}
$
(\mathrm{pr}^{\equiv^{\llbracket2\rrbracket}}_{\boldsymbol{\mathcal{B}}^{(2)}} \circ \mathrm{ip}^{(2,Y)@}_{\boldsymbol{\mathcal{B}}^{(2)}} \circ \mathrm{CH}^{(2)}_{\boldsymbol{\mathcal{B}}^{(2)}})^{\natural}_{\varphi(s)} \left(
\mathrm{sc}_{s}^{0\llbracket\mathbf{Pth}_{\boldsymbol{\mathcal{B}}^{(2)}}^{\mathbf{f}^{(2)}}\rrbracket}\left(
\left\llbracket
\mathfrak{P}^{(2)}
\right\rrbracket_{\varphi(s)}
\right)
\right)
$
\allowdisplaybreaks
\begin{align*}
&=
(\mathrm{pr}^{\equiv^{\llbracket2\rrbracket}}_{\boldsymbol{\mathcal{B}}^{(2)}} \circ \mathrm{ip}^{(2,Y)@}_{\boldsymbol{\mathcal{B}}^{(2)}} \circ \mathrm{CH}^{(2)}_{\boldsymbol{\mathcal{B}}^{(2)}})^{\natural}_{\varphi(s)} \left(
\mathrm{sc}_{\varphi(s)}^{0\llbracket\mathbf{Pth}_{\boldsymbol{\mathcal{B}}^{(2)}}\rrbracket}\left(
\left\llbracket
\mathfrak{P}^{(2)}
\right\rrbracket_{\varphi(s)}
\right)
\right)
\tag{1}
\\
&=
(\mathrm{pr}^{\equiv^{\llbracket2\rrbracket}}_{\boldsymbol{\mathcal{B}}^{(2)}} \circ \mathrm{ip}^{(2,Y)@}_{\boldsymbol{\mathcal{B}}^{(2)}} \circ \mathrm{CH}^{(2)}_{\boldsymbol{\mathcal{B}}^{(2)}})^{\natural}_{\varphi(s)} \left(
\left\llbracket
\mathrm{sc}_{\varphi(s)}^{0\mathbf{Pth}_{\boldsymbol{\mathcal{B}}^{(2)}}}\left(
\mathfrak{P}^{(2)}
\right)
\right\rrbracket_{\varphi(s)}
\right)
\tag{2}
\\
&=
(\mathrm{pr}^{\equiv^{\llbracket2\rrbracket}}_{\boldsymbol{\mathcal{B}}^{(2)}} \circ \mathrm{ip}^{(2,Y)@}_{\boldsymbol{\mathcal{B}}^{(2)}} \circ \mathrm{CH}^{(2)}_{\boldsymbol{\mathcal{B}}^{(2)}})_{\varphi(s)} \left(
\mathrm{sc}_{\varphi(s)}^{0\mathbf{Pth}_{\boldsymbol{\mathcal{B}}^{(2)}}}\left(
\mathfrak{P}^{(2)}
\right)
\right)
\tag{3}
\\
&=
\mathrm{sc}_{\varphi(s)}^{0\mathbf{T}_{\boldsymbol{\mathcal{E}}^{\boldsymbol{\mathcal{B}}^{(2)}}}^{\mathbf{f}^{(2)}}(\mathbf{Pth}_{\boldsymbol{\mathcal{B}}^{(2)}})}\left(
(\mathrm{pr}^{\equiv^{\llbracket2\rrbracket}}_{\boldsymbol{\mathcal{B}}^{(2)}} \circ \mathrm{ip}^{(2,Y)@}_{\boldsymbol{\mathcal{B}}^{(2)}} \circ \mathrm{CH}^{(2)}_{\boldsymbol{\mathcal{B}}^{(2)}})_{\varphi(s)} \left(
\mathfrak{P}^{(2)}
\right)
\right)
\tag{4}
\\
&=
\mathrm{sc}_{\varphi(s)}^{0\mathbf{T}_{\boldsymbol{\mathcal{E}}^{\boldsymbol{\mathcal{B}}^{(2)}}}^{\mathbf{f}^{(2)}}(\mathbf{Pth}_{\boldsymbol{\mathcal{B}}^{(2)}})}\left(
(\mathrm{pr}^{\equiv^{\llbracket2\rrbracket}}_{\boldsymbol{\mathcal{B}}^{(2)}} \circ \mathrm{ip}^{(2,Y)@}_{\boldsymbol{\mathcal{B}}^{(2)}} \circ \mathrm{CH}^{(2)}_{\boldsymbol{\mathcal{B}}^{(2)}})^{\natural}_{\varphi(s)} \left(
\left\llbracket
\mathfrak{P}^{(2)}
\right\rrbracket_{\varphi(s)}
\right)
\right)
\tag{5}
\\
&=
\mathrm{sc}_{s}^{0\mathbf{T}_{\boldsymbol{\mathcal{E}}^{\boldsymbol{\mathcal{B}}^{(2)}}}^{\mathbf{f}^{(2)}}(\mathbf{Pth}_{\boldsymbol{\mathcal{B}}^{(2)}})}\left(
(\mathrm{pr}^{\equiv^{\llbracket2\rrbracket}}_{\boldsymbol{\mathcal{B}}^{(2)}} \circ \mathrm{ip}^{(2,Y)@}_{\boldsymbol{\mathcal{B}}^{(2)}} \circ \mathrm{CH}^{(2)}_{\boldsymbol{\mathcal{B}}^{(2)}})^{\natural}_{\varphi(s)} \left(
\left\llbracket
\mathfrak{P}^{(2)}
\right\rrbracket_{\varphi(s)}
\right)
\right)
\tag{6}
\end{align*}
\end{flushleft}

The first equality unravels the interpretation of the operation symbol $\mathrm{sc}_{s}^{0}$ in the partial $\Sigma^{\boldsymbol{\mathcal{A}}^{(1)}}$-algebra $\llbracket\mathbf{Pth}_{\boldsymbol{\mathcal{B}}^{(2)}}^{\mathbf{f}^{(2)}}\rrbracket$, introduced in Proposition~\ref{PDQPthBDCatAlg};
the second equality unravels the interpretation of the operation symbol $\mathrm{sc}_{s}^{0}$ in the partial $\Sigma^{\boldsymbol{\mathcal{A}}^{(1)}}$-algebra $\llbracket\mathbf{Pth}_{\boldsymbol{\mathcal{B}}^{(2)}}\rrbracket$, introduced in Proposition~\ref{PDVDCatAlg};
the third equality follows from Remark~\ref{RDQPUniv};
the fourth follows from the fact that, according to Proposition~\ref{PDVarKer}, the mapping $(\mathrm{pr}^{\equiv^{\llbracket2\rrbracket}}_{\boldsymbol{\mathcal{B}}^{(2)}} \circ \mathrm{ip}^{(2,Y)@}_{\boldsymbol{\mathcal{B}}^{(2)}} \circ \mathrm{CH}^{(2)}_{\boldsymbol{\mathcal{B}}^{(2)}})_{\varphi(s)}$ is a $\Lambda^{\boldsymbol{\mathcal{B}}^{(2)}}$-homomorphism;
the fifth equality follows from Remark~	\ref{RDQPUniv};
finally, the last equality recovers the interpretation of the operation symbol $\mathrm{sc}_{s}^{0}$ in the partial $\Sigma^{\boldsymbol{\mathcal{A}}^{(1)}}$-algebra $\mathbf{T}_{\boldsymbol{\mathcal{E}}^{\boldsymbol{\mathcal{B}}^{(2)}}}^{\mathbf{f}^{(2)}}(\mathbf{Pth}_{\boldsymbol{\mathcal{B}}^{(2)}})$, introduced in Proposition~\ref{PDFreeBDCatAlg}.

Hence, $(\mathrm{pr}^{\equiv^{\llbracket2\rrbracket}}_{\boldsymbol{\mathcal{B}}^{(2)}} \circ \mathrm{ip}^{(2,Y)@}_{\boldsymbol{\mathcal{B}}^{(2)}} \circ \mathrm{CH}^{(2)}_{\boldsymbol{\mathcal{B}}^{(2)}})^{\natural}_{\varphi}$ is compatible with the $0$-source operation.

{\sffamily The mapping $(\mathrm{pr}^{\equiv^{\llbracket2\rrbracket}}_{\boldsymbol{\mathcal{B}}^{(2)}} \circ \mathrm{ip}^{(2,Y)@}_{\boldsymbol{\mathcal{B}}^{(2)}} \circ \mathrm{CH}^{(2)}_{\boldsymbol{\mathcal{B}}^{(2)}})^{\natural}_{\varphi}$ is compatible with the $0$-target.}

Let $s$ be a sort in $S$ and let us consider the $0$-target operation symbol $\mathrm{tg}_{s}^{0}$ in $\Sigma^{\boldsymbol{\mathcal{A}}^{(1)}}_{s,s}$. Let $\mathfrak{P}^{(2)}$ be a second-order path in $\mathrm{Pth}_{\boldsymbol{\mathcal{B}}^{(2)}, \varphi(s)}$, then the following equality holds
\begin{multline*}
(\mathrm{pr}^{\equiv^{\llbracket2\rrbracket}}_{\boldsymbol{\mathcal{B}}^{(2)}} \circ \mathrm{ip}^{(2,Y)@}_{\boldsymbol{\mathcal{B}}^{(2)}} \circ \mathrm{CH}^{(2)}_{\boldsymbol{\mathcal{B}}^{(2)}})^{\natural}_{\varphi(s)} \left(
\mathrm{tg}_{s}^{0\llbracket\mathbf{Pth}_{\boldsymbol{\mathcal{B}}^{(2)}}^{\mathbf{f}^{(2)}}\rrbracket}\left(
\left\llbracket
\mathfrak{P}^{(2)}
\right\rrbracket_{\varphi(s)}
\right)
\right)
\\
=
\mathrm{tg}_{\varphi(s)}^{0\mathbf{T}_{\boldsymbol{\mathcal{E}}^{\boldsymbol{\mathcal{B}}^{(2)}}}^{\mathbf{f}^{(2)}}(\mathbf{Pth}_{\boldsymbol{\mathcal{B}}^{(2)}})}\left(
(\mathrm{pr}^{\equiv^{\llbracket2\rrbracket}}_{\boldsymbol{\mathcal{B}}^{(2)}} \circ \mathrm{ip}^{(2,Y)@}_{\boldsymbol{\mathcal{B}}^{(2)}} \circ \mathrm{CH}^{(2)}_{\boldsymbol{\mathcal{B}}^{(2)}})^{\natural}_{\varphi(s)} \left(
\left\llbracket
\mathfrak{P}^{(2)}
\right\rrbracket_{\varphi(s)}
\right)
\right).
\end{multline*}

The proof of this case is identical to that of the $0$-source.

Hence, $(\mathrm{pr}^{\equiv^{\llbracket2\rrbracket}}_{\boldsymbol{\mathcal{B}}^{(2)}} \circ \mathrm{ip}^{(2,Y)@}_{\boldsymbol{\mathcal{B}}^{(2)}} \circ \mathrm{CH}^{(2)}_{\boldsymbol{\mathcal{B}}^{(2)}})^{\natural}_{\varphi}$ is compatible with the $0$-target operation.

{\sffamily The mapping $(\mathrm{pr}^{\equiv^{\llbracket2\rrbracket}}_{\boldsymbol{\mathcal{B}}^{(2)}} \circ \mathrm{ip}^{(2,Y)@}_{\boldsymbol{\mathcal{B}}^{(2)}} \circ \mathrm{CH}^{(2)}_{\boldsymbol{\mathcal{B}}^{(2)}})^{\natural}_{\varphi}$ is compatible with the $0$-composition.}

Let $s$ be a sort in $S$ and let us consider the $0$-composition operation symbol $\circ_{s}^{0}$ in $\Sigma_{ss,s}^{\boldsymbol{\mathcal{A}}^{(1)}}$. Let $\mathfrak{P}^{(2)}$ and $\mathfrak{Q}^{(2)}$ be two second-order paths in $\mathrm{Pth}_{\boldsymbol{\mathcal{B}}^{(2)}, \varphi(s)}$ such that
$$
\mathrm{sc}_{\boldsymbol{\mathcal{B}}^{(2)}, \varphi(s)}^{(0,2)}\left(\mathfrak{Q}^{(2)}\right)
=
\mathrm{tg}_{\boldsymbol{\mathcal{B}}^{(2)}, \varphi(s)}^{(0,2)}\left(\mathfrak{P}^{(2)}\right).
$$

Then the following equality holds
\begin{align*}
&
(\mathrm{pr}^{\equiv^{\llbracket2\rrbracket}}_{\boldsymbol{\mathcal{B}}^{(2)}} \circ \mathrm{ip}^{(2,Y)@}_{\boldsymbol{\mathcal{B}}^{(2)}} \circ \mathrm{CH}^{(2)}_{\boldsymbol{\mathcal{B}}^{(2)}})^{\natural}_{\varphi(s)} \left(
\left\llbracket
\mathfrak{Q}^{(2)}
\right\rrbracket_{\varphi(s)}
\circ_{s}^{0\llbracket\mathbf{Pth}_{\boldsymbol{\mathcal{B}}^{(2)}}^{\mathbf{f}^{(2)}}\rrbracket}
\left\llbracket
\mathfrak{P}^{(2)}
\right\rrbracket_{\varphi(s)}
\right)
\\
&\hspace{1.5cm}
=(\mathrm{pr}^{\equiv^{\llbracket2\rrbracket}}_{\boldsymbol{\mathcal{B}}^{(2)}} \circ \mathrm{ip}^{(2,Y)@}_{\boldsymbol{\mathcal{B}}^{(2)}} \circ \mathrm{CH}^{(2)}_{\boldsymbol{\mathcal{B}}^{(2)}})^{\natural}_{\varphi(s)} \left(
\left\llbracket
\mathfrak{Q}^{(2)}
\right\rrbracket_{\varphi(s)}
\right)
\\
&\hspace{2cm}
\circ_{s}^{0\mathbf{T}_{\boldsymbol{\mathcal{E}}^{\boldsymbol{\mathcal{B}}^{(2)}}}^{\mathbf{f}^{(2)}}(\mathbf{Pth}_{\boldsymbol{\mathcal{B}}^{(2)}})}
(\mathrm{pr}^{\equiv^{\llbracket2\rrbracket}}_{\boldsymbol{\mathcal{B}}^{(2)}} \circ \mathrm{ip}^{(2,Y)@}_{\boldsymbol{\mathcal{B}}^{(2)}} \circ \mathrm{CH}^{(2)}_{\boldsymbol{\mathcal{B}}^{(2)}})^{\natural}_{\varphi(s)} \left(
\left\llbracket
\mathfrak{P}^{(2)}
\right\rrbracket_{\varphi(s)}
\right).
\end{align*}

The proof of this case is identical to that of the $0$-source.

Hence, $(\mathrm{pr}^{\equiv^{\llbracket2\rrbracket}}_{\boldsymbol{\mathcal{B}}^{(2)}} \circ \mathrm{ip}^{(2,Y)@}_{\boldsymbol{\mathcal{B}}^{(2)}} \circ \mathrm{CH}^{(2)}_{\boldsymbol{\mathcal{B}}^{(2)}})^{\natural}_{\varphi}$ is compatible with the $0$-composition operation.

{\sffamily The mapping $(\mathrm{pr}^{\equiv^{\llbracket2\rrbracket}}_{\boldsymbol{\mathcal{B}}^{(2)}} \circ \mathrm{ip}^{(2,Y)@}_{\boldsymbol{\mathcal{B}}^{(2)}} \circ \mathrm{CH}^{(2)}_{\boldsymbol{\mathcal{B}}^{(2)}})^{\natural}_{\varphi}$ is compatible with the second-order rewrite rules.}

Let $s$ be a sort in $S$ and $\mathfrak{p}^{(2)}$ a rewrite rule in $\mathcal{A}_{s}^{(2)}$. Thus, the following chain of equalities holds
\begin{flushleft}
$
(\mathrm{pr}^{\equiv^{\llbracket2\rrbracket}}_{\boldsymbol{\mathcal{B}}^{(2)}} \circ \mathrm{ip}^{(2,Y)@}_{\boldsymbol{\mathcal{B}}^{(2)}} \circ \mathrm{CH}^{(2)}_{\boldsymbol{\mathcal{B}}^{(2)}})^{\natural}_{\varphi(s)} \left(
\mathfrak{p}^{(2)\llbracket\mathbf{Pth}_{\boldsymbol{\mathcal{B}}^{(2)}}^{\mathbf{f}^{(2)}}\rrbracket}
\right)
$
\begin{align*}
&=
(\mathrm{pr}^{\equiv^{\llbracket2\rrbracket}}_{\boldsymbol{\mathcal{B}}^{(2)}} \circ \mathrm{ip}^{(2,Y)@}_{\boldsymbol{\mathcal{B}}^{(2)}} \circ \mathrm{CH}^{(2)}_{\boldsymbol{\mathcal{B}}^{(2)}})^{\natural}_{\varphi(s)} \left(
\left\llbracket
f^{(2)}\left(
\mathfrak{p}^{(2)}
\right)
\right\rrbracket_{\varphi(s)}
\right)
\tag{1}
\\
&=
(\mathrm{pr}^{\equiv^{\llbracket2\rrbracket}}_{\boldsymbol{\mathcal{B}}^{(2)}} \circ \mathrm{ip}^{(2,Y)@}_{\boldsymbol{\mathcal{B}}^{(2)}} \circ \mathrm{CH}^{(2)}_{\boldsymbol{\mathcal{B}}^{(2)}})_{\varphi(s)} \left(
f^{(2)}\left(
\mathfrak{p}^{(2)}
\right)
\right)
\tag{2}
\\
&=
\mathrm{pr}^{\equiv^{\llbracket2\rrbracket}}_{\boldsymbol{\mathcal{B}}^{(2)}, \varphi(s)} 
\circ
\mathrm{ip}^{(2,Y)@}_{\boldsymbol{\mathcal{B}}^{(2)}, \varphi(s)}
\circ
\mathrm{CH}^{(2)}_{\boldsymbol{\mathcal{B}}^{(2)}, \varphi(s)} \left(
f^{(2)}\left(
\mathfrak{p}^{(2)}
\right)
\right)
\tag{3}
\\
&=
\mathfrak{p}^{(2)\mathbf{T}_{\boldsymbol{\mathcal{E}}^{\boldsymbol{\mathcal{B}}^{(2)}}}^{\mathbf{f}^{(2)}}(\mathbf{Pth}_{\boldsymbol{\mathcal{B}}^{(2)}})}.
\tag{4}
\end{align*}
\end{flushleft}

The first equality unravels the interpretation of the constant $\mathfrak{p}^{(2)}$ in the partial $\Sigma^{\boldsymbol{\mathcal{A}}^{(2)}}$-algebra $\llbracket\mathbf{Pth}^{\mathbf{f}^{(2)}}_{\boldsymbol{\mathcal{B}}^{(2)}}\rrbracket$, introduced in Proposition~\ref{PDQPthBDCatAlg};
the second equality follows from Remark~\ref{RDQPUniv};
the third equality unravels the definition of the $\varphi(s)$ component of a composition of many-sorted mappings;
finally, the last equality recovers the interpretation of the constant $\mathfrak{p}^{(2)}$ in the partial $\Sigma^{\boldsymbol{\mathcal{A}}^{(2)}}$-algebra $\mathbf{T}_{\boldsymbol{\mathcal{E}}^{\boldsymbol{\mathcal{B}}^{(2)}}}^{\mathbf{f}^{(2)}}(\mathbf{Pth}_{\boldsymbol{\mathcal{B}}^{(2)}})$, introduced in Proposition~\ref{PDFreeBDCatAlg}.

Hence, $(\mathrm{pr}^{\equiv^{\llbracket2\rrbracket}}_{\boldsymbol{\mathcal{B}}^{(2)}} \circ \mathrm{ip}^{(2,Y)@}_{\boldsymbol{\mathcal{B}}^{(2)}} \circ \mathrm{CH}^{(2)}_{\boldsymbol{\mathcal{B}}^{(2)}})^{\natural}_{\varphi}$ is compatible with the second-order rewrite rules.

{\sffamily The mapping $(\mathrm{pr}^{\equiv^{\llbracket2\rrbracket}}_{\boldsymbol{\mathcal{B}}^{(2)}} \circ \mathrm{ip}^{(2,Y)@}_{\boldsymbol{\mathcal{B}}^{(2)}} \circ \mathrm{CH}^{(2)}_{\boldsymbol{\mathcal{B}}^{(2)}})^{\natural}_{\varphi}$ is compatible with the $1$-source.}

Let $s$ be a sort in $S$ and let us consider the $1$-source operation symbol $\mathrm{sc}_{s}^{1}$ in $\Sigma^{\boldsymbol{\mathcal{A}}^{(2)}}_{s,s}$. Let $\mathfrak{P}^{(2)}$ be a second-order path in $\mathrm{Pth}_{\boldsymbol{\mathcal{B}}^{(2)}, \varphi(s)}$.

The following chain of equalities holds
\begin{flushleft}
$
(\mathrm{pr}^{\equiv^{\llbracket2\rrbracket}}_{\boldsymbol{\mathcal{B}}^{(2)}} \circ \mathrm{ip}^{(2,Y)@}_{\boldsymbol{\mathcal{B}}^{(2)}} \circ \mathrm{CH}^{(2)}_{\boldsymbol{\mathcal{B}}^{(2)}})^{\natural}_{\varphi(s)} \left(
\mathrm{sc}_{s}^{1\llbracket\mathbf{Pth}_{\boldsymbol{\mathcal{B}}^{(2)}}^{\mathbf{f}^{(2)}}\rrbracket}\left(
\left\llbracket
\mathfrak{P}^{(2)}
\right\rrbracket_{\varphi(s)}
\right)
\right)
$
\allowdisplaybreaks
\begin{align*}
&=
(\mathrm{pr}^{\equiv^{\llbracket2\rrbracket}}_{\boldsymbol{\mathcal{B}}^{(2)}} \circ \mathrm{ip}^{(2,Y)@}_{\boldsymbol{\mathcal{B}}^{(2)}} \circ \mathrm{CH}^{(2)}_{\boldsymbol{\mathcal{B}}^{(2)}})^{\natural}_{\varphi(s)} \left(
\mathrm{sc}_{\varphi(s)}^{1\llbracket\mathbf{Pth}_{\boldsymbol{\mathcal{B}}^{(2)}}\rrbracket}\left(
\left\llbracket
\mathfrak{P}^{(2)}
\right\rrbracket_{\varphi(s)}
\right)
\right)
\tag{1}
\\
&=
(\mathrm{pr}^{\equiv^{\llbracket2\rrbracket}}_{\boldsymbol{\mathcal{B}}^{(2)}} \circ \mathrm{ip}^{(2,Y)@}_{\boldsymbol{\mathcal{B}}^{(2)}} \circ \mathrm{CH}^{(2)}_{\boldsymbol{\mathcal{B}}^{(2)}})^{\natural}_{\varphi(s)} \left(
\left\llbracket
\mathrm{sc}_{\varphi(s)}^{1\mathbf{Pth}_{\boldsymbol{\mathcal{B}}^{(2)}}}\left(
\mathfrak{P}^{(2)}
\right)
\right\rrbracket_{\varphi(s)}
\right)
\tag{2}
\\
&=
(\mathrm{pr}^{\equiv^{\llbracket2\rrbracket}}_{\boldsymbol{\mathcal{B}}^{(2)}} \circ \mathrm{ip}^{(2,Y)@}_{\boldsymbol{\mathcal{B}}^{(2)}} \circ \mathrm{CH}^{(2)}_{\boldsymbol{\mathcal{B}}^{(2)}})_{\varphi(s)} \left(
\mathrm{sc}_{\varphi(s)}^{1\mathbf{Pth}_{\boldsymbol{\mathcal{B}}^{(2)}}}\left(
\mathfrak{P}^{(2)}
\right)
\right)
\tag{3}
\\
&=
\mathrm{sc}_{\varphi(s)}^{1\mathbf{T}_{\boldsymbol{\mathcal{E}}^{\boldsymbol{\mathcal{B}}^{(2)}}}^{\mathbf{f}^{(2)}}(\mathbf{Pth}_{\boldsymbol{\mathcal{B}}^{(2)}})}\left(
(\mathrm{pr}^{\equiv^{\llbracket2\rrbracket}}_{\boldsymbol{\mathcal{B}}^{(2)}} \circ \mathrm{ip}^{(2,Y)@}_{\boldsymbol{\mathcal{B}}^{(2)}} \circ \mathrm{CH}^{(2)}_{\boldsymbol{\mathcal{B}}^{(2)}})_{\varphi(s)} \left(
\mathfrak{P}^{(2)}
\right)
\right)
\tag{4}
\\
&=
\mathrm{sc}_{\varphi(s)}^{1\mathbf{T}_{\boldsymbol{\mathcal{E}}^{\boldsymbol{\mathcal{B}}^{(2)}}}^{\mathbf{f}^{(2)}}(\mathbf{Pth}_{\boldsymbol{\mathcal{B}}^{(2)}})}\left(
(\mathrm{pr}^{\equiv^{\llbracket2\rrbracket}}_{\boldsymbol{\mathcal{B}}^{(2)}} \circ \mathrm{ip}^{(2,Y)@}_{\boldsymbol{\mathcal{B}}^{(2)}} \circ \mathrm{CH}^{(2)}_{\boldsymbol{\mathcal{B}}^{(2)}})^{\natural}_{\varphi(s)} \left(
\left\llbracket
\mathfrak{P}^{(2)}
\right\rrbracket_{\varphi(s)}
\right)
\right)
\tag{5}
\\
&=
\mathrm{sc}_{s}^{1\mathbf{T}_{\boldsymbol{\mathcal{E}}^{\boldsymbol{\mathcal{B}}^{(2)}}}^{\mathbf{f}^{(2)}}(\mathbf{Pth}_{\boldsymbol{\mathcal{B}}^{(2)}})}\left(
(\mathrm{pr}^{\equiv^{\llbracket2\rrbracket}}_{\boldsymbol{\mathcal{B}}^{(2)}} \circ \mathrm{ip}^{(2,Y)@}_{\boldsymbol{\mathcal{B}}^{(2)}} \circ \mathrm{CH}^{(2)}_{\boldsymbol{\mathcal{B}}^{(2)}})^{\natural}_{\varphi(s)} \left(
\left\llbracket
\mathfrak{P}^{(2)}
\right\rrbracket_{\varphi(s)}
\right)
\right)
\tag{6}
\end{align*}
\end{flushleft}

The first equality unravels the interpretation of the operation symbol $\mathrm{sc}_{s}^{1}$ in the partial $\Sigma^{\boldsymbol{\mathcal{A}}^{(1)}}$-algebra $\llbracket\mathbf{Pth}_{\boldsymbol{\mathcal{B}}^{(2)}}^{\mathbf{f}^{(2)}}\rrbracket$, introduced in Proposition~\ref{PDQPthBDCatAlg};
the second equality unravels the interpretation of the operation symbol $\mathrm{sc}_{s}^{1}$ in the partial $\Sigma^{\boldsymbol{\mathcal{A}}^{(1)}}$-algebra $\llbracket\mathbf{Pth}_{\boldsymbol{\mathcal{B}}^{(2)}}\rrbracket$, introduced in Proposition~\ref{PDVDCatAlg};
the third equality follows from Remark~\ref{RDQPUniv};
the fourth follows from the fact that, according to Proposition~\ref{PDVarKer}, the mapping $(\mathrm{pr}^{\equiv^{\llbracket2\rrbracket}}_{\boldsymbol{\mathcal{B}}^{(2)}} \circ \mathrm{ip}^{(2,Y)@}_{\boldsymbol{\mathcal{B}}^{(2)}} \circ \mathrm{CH}^{(2)}_{\boldsymbol{\mathcal{B}}^{(2)}})_{\varphi(s)}$ is a $\Lambda^{\boldsymbol{\mathcal{B}}^{(2)}}$-homomorphism;
the fifth equality follows from Remark~	\ref{RDQPUniv};
finally, the last equality recovers the interpretation of the operation symbol $\mathrm{sc}_{s}^{1}$ in the partial $\Sigma^{\boldsymbol{\mathcal{A}}^{(1)}}$-algebra $\mathbf{T}_{\boldsymbol{\mathcal{E}}^{\boldsymbol{\mathcal{B}}^{(2)}}}^{\mathbf{f}^{(2)}}(\mathbf{Pth}_{\boldsymbol{\mathcal{B}}^{(2)}})$, introduced in Proposition~\ref{PDFreeBDCatAlg}.

Hence, $(\mathrm{pr}^{\equiv^{\llbracket2\rrbracket}}_{\boldsymbol{\mathcal{B}}^{(2)}} \circ \mathrm{ip}^{(2,Y)@}_{\boldsymbol{\mathcal{B}}^{(2)}} \circ \mathrm{CH}^{(2)}_{\boldsymbol{\mathcal{B}}^{(2)}})^{\natural}_{\varphi}$ is compatible with the $1$-source operation.

{\sffamily The mapping $(\mathrm{pr}^{\equiv^{\llbracket2\rrbracket}}_{\boldsymbol{\mathcal{B}}^{(2)}} \circ \mathrm{ip}^{(2,Y)@}_{\boldsymbol{\mathcal{B}}^{(2)}} \circ \mathrm{CH}^{(2)}_{\boldsymbol{\mathcal{B}}^{(2)}})^{\natural}_{\varphi}$ is compatible with the $1$-target.}

Let $s$ be a sort in $S$ and let us consider the $1$-target operation symbol $\mathrm{tg}_{s}^{1}$ in $\Sigma^{\boldsymbol{\mathcal{A}}^{(2)}}_{s,s}$. Let $\mathfrak{P}^{(2)}$ be a second-order path in $\mathrm{Pth}_{\boldsymbol{\mathcal{B}}^{(2)}, \varphi(s)}$, then the following equality holds
\begin{multline*}
(\mathrm{pr}^{\equiv^{\llbracket2\rrbracket}}_{\boldsymbol{\mathcal{B}}^{(2)}} \circ \mathrm{ip}^{(2,Y)@}_{\boldsymbol{\mathcal{B}}^{(2)}} \circ \mathrm{CH}^{(2)}_{\boldsymbol{\mathcal{B}}^{(2)}})^{\natural}_{\varphi(s)} \left(
\mathrm{tg}_{s}^{1\llbracket\mathbf{Pth}_{\boldsymbol{\mathcal{B}}^{(2)}}^{\mathbf{f}^{(2)}}\rrbracket}\left(
\left\llbracket
\mathfrak{P}^{(2)}
\right\rrbracket_{\varphi(s)}
\right)
\right)
\\
=
\mathrm{tg}_{\varphi(s)}^{1\mathbf{T}_{\boldsymbol{\mathcal{E}}^{\boldsymbol{\mathcal{B}}^{(2)}}}^{\mathbf{f}^{(2)}}(\mathbf{Pth}_{\boldsymbol{\mathcal{B}}^{(2)}})}\left(
(\mathrm{pr}^{\equiv^{\llbracket2\rrbracket}}_{\boldsymbol{\mathcal{B}}^{(2)}} \circ \mathrm{ip}^{(2,Y)@}_{\boldsymbol{\mathcal{B}}^{(2)}} \circ \mathrm{CH}^{(2)}_{\boldsymbol{\mathcal{B}}^{(2)}})^{\natural}_{\varphi(s)} \left(
\left\llbracket
\mathfrak{P}^{(2)}
\right\rrbracket_{\varphi(s)}
\right)
\right).
\end{multline*}

The proof of this case is identical to that of the $1$-source.

Hence, $(\mathrm{pr}^{\equiv^{\llbracket2\rrbracket}}_{\boldsymbol{\mathcal{B}}^{(2)}} \circ \mathrm{ip}^{(2,Y)@}_{\boldsymbol{\mathcal{B}}^{(2)}} \circ \mathrm{CH}^{(2)}_{\boldsymbol{\mathcal{B}}^{(2)}})^{\natural}_{\varphi}$ is compatible with the $1$-target operation.

{\sffamily The mapping $(\mathrm{pr}^{\equiv^{\llbracket2\rrbracket}}_{\boldsymbol{\mathcal{B}}^{(2)}} \circ \mathrm{ip}^{(2,Y)@}_{\boldsymbol{\mathcal{B}}^{(2)}} \circ \mathrm{CH}^{(2)}_{\boldsymbol{\mathcal{B}}^{(2)}})^{\natural}_{\varphi}$ is compatible with the $1$-composition.}

Let $s$ be a sort in $S$ and let us consider the $1$-composition operation symbol $\circ_{s}^{1}$ in $\Sigma_{ss,s}^{\boldsymbol{\mathcal{A}}^{(2)}}$. Let $\mathfrak{P}^{(2)}$ and $\mathfrak{Q}^{(2)}$ be two second-order paths in $\mathrm{Pth}_{\boldsymbol{\mathcal{B}}^{(2)}, \varphi(s)}$ such that
$$
\mathrm{sc}_{\boldsymbol{\mathcal{B}}^{(2)}, \varphi(s)}^{([1],2)}\left(\mathfrak{Q}^{(2)}\right)
=
\mathrm{tg}_{\boldsymbol{\mathcal{B}}^{(2)}, \varphi(s)}^{([1],2)}\left(\mathfrak{P}^{(2)}\right).
$$
Then the following equality holds
\begin{align*}
&
(\mathrm{pr}^{\equiv^{\llbracket2\rrbracket}}_{\boldsymbol{\mathcal{B}}^{(2)}} \circ \mathrm{ip}^{(2,Y)@}_{\boldsymbol{\mathcal{B}}^{(2)}} \circ \mathrm{CH}^{(2)}_{\boldsymbol{\mathcal{B}}^{(2)}})^{\natural}_{\varphi(s)} \left(
\left\llbracket
\mathfrak{Q}^{(2)}
\right\rrbracket_{\varphi(s)}
\circ_{s}^{1\llbracket\mathbf{Pth}_{\boldsymbol{\mathcal{B}}^{(2)}}^{\mathbf{f}^{(2)}}\rrbracket}
\left\llbracket
\mathfrak{P}^{(2)}
\right\rrbracket_{\varphi(s)}
\right)
\\
&\hspace{1.5cm}
=(\mathrm{pr}^{\equiv^{\llbracket2\rrbracket}}_{\boldsymbol{\mathcal{B}}^{(2)}} \circ \mathrm{ip}^{(2,Y)@}_{\boldsymbol{\mathcal{B}}^{(2)}} \circ \mathrm{CH}^{(2)}_{\boldsymbol{\mathcal{B}}^{(2)}})^{\natural}_{\varphi(s)} \left(
\left\llbracket
\mathfrak{Q}^{(2)}
\right\rrbracket_{\varphi(s)}
\right)
\\
&\hspace{2cm}
\circ_{s}^{1\mathbf{T}_{\boldsymbol{\mathcal{E}}^{\boldsymbol{\mathcal{B}}^{(2)}}}^{\mathbf{f}^{(2)}}(\mathbf{Pth}_{\boldsymbol{\mathcal{B}}^{(2)}})}
(\mathrm{pr}^{\equiv^{\llbracket2\rrbracket}}_{\boldsymbol{\mathcal{B}}^{(2)}} \circ \mathrm{ip}^{(2,Y)@}_{\boldsymbol{\mathcal{B}}^{(2)}} \circ \mathrm{CH}^{(2)}_{\boldsymbol{\mathcal{B}}^{(2)}})^{\natural}_{\varphi(s)} \left(
\left\llbracket
\mathfrak{P}^{(2)}
\right\rrbracket_{\varphi(s)}
\right).
\end{align*}

The proof of this case is identical to that of the $1$-source.

Hence, $(\mathrm{pr}^{\equiv^{\llbracket2\rrbracket}}_{\boldsymbol{\mathcal{B}}^{(2)}} \circ \mathrm{ip}^{(2,Y)@}_{\boldsymbol{\mathcal{B}}^{(2)}} \circ \mathrm{CH}^{(2)}_{\boldsymbol{\mathcal{B}}^{(2)}})^{\natural}_{\varphi}$ is compatible with the $1$-composition operation.

This completes the proof.
\end{proof}

\begin{proposition}
\label{PDQIpfcBDCatHom}
Let $\mathbf{f}^{(2)}=(\varphi, c, (f^{(i)})_{i\in 3})$ be a second-order morphism from $\boldsymbol{\mathcal{A}}^{(1)}$ to $\boldsymbol{\mathcal{B}}^{(1)}$. Then the mapping 
$$
\mathrm{pr}^{\llbracket\cdot\rrbracket\mathsf{p}}_{\boldsymbol{\mathcal{B}}^{(2)}, \varphi}
\colon 
\T_{\boldsymbol{\mathcal{E}}^{\boldsymbol{\mathcal{B}}^{(2)}}}(\mathbf{Pth}_{\boldsymbol{\mathcal{B}}^{(2)}})_{\varphi}
\mor
\llbracket\mathrm{Pth}_{\boldsymbol{\mathcal{B}}^{(2)}}\rrbracket_{\varphi}
$$
is a $\Sigma^{\boldsymbol{\mathcal{A}}^{(2)}}$-homomorphism from $\mathbf{T}_{\boldsymbol{\mathcal{E}}^{\boldsymbol{\mathcal{B}}^{(2)}}}^{\mathbf{f}^{(2)}}(\mathbf{Pth}_{\boldsymbol{\mathcal{B}}^{(2)}})$ to $\llbracket\mathbf{Pth}_{\boldsymbol{\mathcal{B}}^{(2)}}^{\mathbf{f}^{(2)}}\rrbracket$.
\end{proposition}

\begin{proof}
We prove that $\mathrm{pr}^{\llbracket\cdot\rrbracket\mathsf{p}}_{\boldsymbol{\mathcal{B}}^{(2)}, \varphi}$ is compatible with every operation symbol in $\Sigma^{\boldsymbol{\mathcal{A}}^{(2)}}$.

{\sffamily The mapping $\mathrm{pr}^{\llbracket\cdot\rrbracket\mathsf{p}}_{\boldsymbol{\mathcal{B}}^{(2)}, \varphi}$ is a $\Sigma$-homomorphism.}

Note that $\mathrm{pr}^{\llbracket\cdot\rrbracket\mathsf{p}}_{\boldsymbol{\mathcal{B}}^{(2)}, \varphi} = \mathbf{c}_{\mathfrak{d}}^{\ast} (\mathrm{pr}^{\llbracket\cdot\rrbracket\mathsf{p}}_{\boldsymbol{\mathcal{B}}^{(2)}})$. By Corollary~\ref{CDVarPr}, the mapping $\mathrm{pr}^{\llbracket\cdot\rrbracket\mathsf{p}}_{\boldsymbol{\mathcal{B}}^{(2)}}$ is a $\Lambda^{\boldsymbol{\mathcal{B}}^{(2)}}$-homomorphism, thus in particular a $\Lambda$-homomorphism. Therefore, it follows from Proposition~\ref{PFunSig} that the mapping $\mathrm{pr}^{\llbracket\cdot\rrbracket\mathsf{p}}_{\boldsymbol{\mathcal{B}}^{(2)}, \varphi}$ is a $\Sigma$-homomorphism.


{\sffamily The mapping $\mathrm{pr}^{\llbracket\cdot\rrbracket\mathsf{p}}_{\boldsymbol{\mathcal{B}}^{(2)}, \varphi}$ is compatible with the first-order rewrite rules.}

Let $s$ be a sort in $S$ and $\mathfrak{p}$ a rewrite rule in $\mathcal{A}_{s}^{(1)}$. Thus, the following chain of equalities holds
\begin{flushleft}
$
\mathrm{pr}^{\llbracket\cdot\rrbracket\mathsf{p}}_{\boldsymbol{\mathcal{B}}^{(2)}, \varphi(s)} \left(
\mathfrak{p}^{\mathbf{T}_{\boldsymbol{\mathcal{E}}^{\boldsymbol{\mathcal{B}}^{(2)}}}^{\mathbf{f}^{(2)}}(\mathbf{Pth}_{\boldsymbol{\mathcal{B}}^{(2)}})}
\right)
$
\allowdisplaybreaks
\begin{align*}
&=
\mathrm{pr}^{\llbracket\cdot\rrbracket\mathsf{p}}_{\boldsymbol{\mathcal{B}}^{(2)}, \varphi(s)} \left(
\mathrm{pr}^{\equiv^{\llbracket2\rrbracket}}_{\boldsymbol{\mathcal{B}}^{(2)}, \varphi(s)} \circ \mathrm{ip}^{(2,Y)@}_{\boldsymbol{\mathcal{B}}^{(2)}, \varphi(s)} \circ \mathrm{CH}^{(2)}_{\boldsymbol{\mathcal{B}}^{(2)}, \varphi(s)} \left(
f^{(2)\flat}_{s}\left(
\mathfrak{p}^{\mathbf{Pth}_{\boldsymbol{\mathcal{A}}^{(2)}}}
\right)
\right)
\right)
\tag{1}
\\
&=
\mathrm{pr}^{\llbracket\cdot\rrbracket\mathsf{p}}_{\boldsymbol{\mathcal{B}}^{(2)}, \varphi(s)} \left(
(\mathrm{pr}^{\equiv^{\llbracket2\rrbracket}}_{\boldsymbol{\mathcal{B}}^{(2)}} \circ \mathrm{ip}^{(2,Y)@}_{\boldsymbol{\mathcal{B}}^{(2)}} \circ \mathrm{CH}^{(2)}_{\boldsymbol{\mathcal{B}}^{(2)}})_{\varphi(s)} \left(
f^{(2)\flat}_{s}\left(
\mathfrak{p}^{\mathbf{Pth}_{\boldsymbol{\mathcal{A}}^{(2)}}}
\right)
\right)
\right)
\tag{2}
\\
&=
\resizebox{.89\textwidth}{!}{%
$
\mathrm{pr}^{\llbracket\cdot\rrbracket\mathsf{p}}_{\boldsymbol{\mathcal{B}}^{(2)}, \varphi(s)} \left(
(\mathrm{pr}^{\equiv^{\llbracket2\rrbracket}}_{\boldsymbol{\mathcal{B}}^{(2)}} \circ \mathrm{ip}^{(2,Y)@}_{\boldsymbol{\mathcal{B}}^{(2)}} \circ \mathrm{CH}^{(2)}_{\boldsymbol{\mathcal{B}}^{(2)}})_{\varphi(s)}^{\natural} \left(
\left\llbracket
f^{(2)\flat}_{s}\left(
\mathfrak{p}^{\mathbf{Pth}_{\boldsymbol{\mathcal{A}}^{(2)}}}
\right)
\right\rrbracket_{\varphi(s)}
\right)
\right)
$
}
\tag{3}
\\
&=
\left\llbracket
f^{(2)\flat}_{s}\left(
\mathfrak{p}^{\mathbf{Pth}_{\boldsymbol{\mathcal{A}}^{(2)}}}
\right)
\right\rrbracket_{\varphi(s)}
\tag{4}
\\
&=
\mathfrak{p}^{\llbracket\mathbf{Pth}_{\boldsymbol{\mathcal{B}}^{(2)}}^{\mathbf{f}^{(2)}}\rrbracket}.
\end{align*}
\end{flushleft}

The first equality unravels the interpretation of the constant $\mathfrak{p}$ in the partial $\Sigma^{\boldsymbol{\mathcal{A}}^{(2)}}$-algebra $\mathbf{T}_{\boldsymbol{\mathcal{E}}^{\boldsymbol{\mathcal{B}}^{(2)}}}^{\mathbf{f}^{(2)}}(\mathbf{Pth}_{\boldsymbol{\mathcal{B}}^{(2)}})$, introduced in Proposition~\ref{PDFreeBDCatAlg};
the second equality recovers the definition of the $\varphi(s)$ component of a composition of many-sorted mappings;
the third equality follows from Remark~\ref{RDQPUniv};
The fourth equality follows from Theorem~\ref{TDPthFree};
finally, the last equality recovers the interpretation of the constant $\mathfrak{p}$ in the partial $\Sigma^{\boldsymbol{\mathcal{A}}^{(2)}}$-algebra $\llbracket\mathbf{Pth}_{\boldsymbol{\mathcal{B}}^{(2)}}^{\mathbf{f}^{(2)}}\rrbracket$, introduced in Proposition~\ref{PDQPthBDCatAlg}.

Hence, $\mathrm{pr}^{\llbracket\cdot\rrbracket\mathsf{p}}_{\boldsymbol{\mathcal{B}}^{(2)}, \varphi}$ is compatible with the first-order rewrite rules.

{\sffamily The mapping $\mathrm{pr}^{\llbracket\cdot\rrbracket\mathsf{p}}_{\boldsymbol{\mathcal{B}}^{(2)}, \varphi}$ is compatible with the $0$-source.}

Let $s$ be a sort in $S$ and let us consider the $0$-source operation symbol $\mathrm{sc}_{s}^{0}$ in $\Sigma^{\boldsymbol{\mathcal{A}}^{(2)}}_{s,s}$. Let $x$ be an element in $\T_{\boldsymbol{\mathcal{E}}^{\boldsymbol{\mathcal{B}}^{(2)}}}(\mathbf{Pth}_{\boldsymbol{\mathcal{B}}^{(2)}})_{\varphi(s)}$.

The following chain of equalities holds
\allowdisplaybreaks
\begin{align*}
\mathrm{pr}^{\llbracket\cdot\rrbracket\mathsf{p}}_{\boldsymbol{\mathcal{B}}^{(2)}, \varphi(s)} \left(
\mathrm{sc}_{s}^{0\mathbf{T}_{\boldsymbol{\mathcal{E}}^{\boldsymbol{\mathcal{B}}^{(2)}}}^{\mathbf{f}^{(2)}}(\mathbf{Pth}_{\boldsymbol{\mathcal{B}}^{(2)}})}\left(
x
\right)
\right)
&=
\mathrm{pr}^{\llbracket\cdot\rrbracket\mathsf{p}}_{\boldsymbol{\mathcal{B}}^{(2)}, \varphi(s)} \left(
\mathrm{sc}_{\varphi(s)}^{0\mathbf{T}_{\boldsymbol{\mathcal{E}}^{\boldsymbol{\mathcal{B}}^{(2)}}}(\mathbf{Pth}_{\boldsymbol{\mathcal{B}}^{(2)}})}\left(
x
\right)
\right)
\tag{1}
\\
&=
\mathrm{sc}_{\varphi(s)}^{0\llbracket\mathbf{Pth}_{\boldsymbol{\mathcal{B}}^{(2)}}\rrbracket}\left(
\mathrm{pr}^{\llbracket\cdot\rrbracket\mathsf{p}}_{\boldsymbol{\mathcal{B}}^{(2)}, \varphi(s)} \left(
x
\right)
\right)
\tag{2}
\\
&=
\mathrm{sc}_{s}^{0\llbracket\mathbf{Pth}_{\boldsymbol{\mathcal{B}}^{(2)}}^{\mathbf{f}^{(2)}}\rrbracket}\left(
\mathrm{pr}^{\llbracket\cdot\rrbracket\mathsf{p}}_{\boldsymbol{\mathcal{B}}^{(2)}, \varphi(s)} \left(
x
\right)
\right).
\tag{3}
\end{align*}
The first equality unravels the interpretation of the operation symbol $\mathrm{sc}_{s}^{0}$ in the partial $\Sigma^{\boldsymbol{\mathcal{A}}^{(2)}}$-algebra $\mathbf{T}_{\boldsymbol{\mathcal{E}}^{\boldsymbol{\mathcal{B}}^{(2)}}}^{\mathbf{f}^{(2)}}(\mathbf{Pth}_{\boldsymbol{\mathcal{B}}^{(2)}})$, introduced in Proposition~\ref{PDFreeBDCatAlg};
the second equality follows from the fact that, according to Proposition~\ref{CDVarPr}, $\mathrm{pr}^{\llbracket\cdot\rrbracket\mathsf{p}}_{\boldsymbol{\mathcal{B}}^{(2)}}$ is a $\Lambda^{\boldsymbol{\mathcal{B}}^{(2)}}$-homomorphism from $\mathbf{T}_{\boldsymbol{\mathcal{E}}^{\boldsymbol{\mathcal{B}}^{(2)}}}(\mathbf{Pth}_{\boldsymbol{\mathcal{B}}^{(2)}})$ to $\llbracket\mathbf{Pth}_{\boldsymbol{\mathcal{B}}^{(2)}}\rrbracket$;
finally, the last equality recovers the interpretation of the operation symbol $\mathrm{sc}_{s}^{0}$ in the partial $\Sigma^{\boldsymbol{\mathcal{A}}^{(2)}}$-algebra $\llbracket\mathbf{Pth}_{\boldsymbol{\mathcal{B}}^{(2)}}^{\mathbf{f}^{(2)}}\rrbracket$, introduced in Proposition~\ref{PDQPthBDCatAlg}.

Hence, $\mathrm{pr}^{\llbracket\cdot\rrbracket\mathsf{p}}_{\boldsymbol{\mathcal{B}}^{(2)}, \varphi}$ is compatible with the $0$-source operation.

{\sffamily The mapping $\mathrm{pr}^{\llbracket\cdot\rrbracket\mathsf{p}}_{\boldsymbol{\mathcal{B}}^{(2)}, \varphi}$ is compatible with the $0$-target.}

Let $s$ be a sort in $S$ and let us consider the $0$-target operation symbol $\mathrm{tg}_{s}^{0}$ in $\Sigma^{\boldsymbol{\mathcal{A}}^{(2)}}_{s,s}$. Let $x$ be an element in $\T_{\boldsymbol{\mathcal{E}}^{\boldsymbol{\mathcal{B}}^{(2)}}}(\mathbf{Pth}_{\boldsymbol{\mathcal{B}}^{(2)}})_{\varphi(s)}$, then the following equality holds
$$
\mathrm{pr}^{\llbracket\cdot\rrbracket\mathsf{p}}_{\boldsymbol{\mathcal{B}}^{(2)}, \varphi(s)} \left(
\mathrm{tg}_{s}^{0\mathbf{T}_{\boldsymbol{\mathcal{E}}^{\boldsymbol{\mathcal{B}}^{(2)}}}^{\mathbf{f}^{(2)}}(\mathbf{Pth}_{\boldsymbol{\mathcal{B}}^{(2)}})}\left(
x
\right)
\right)
=
\mathrm{tg}_{s}^{0\llbracket\mathbf{Pth}_{\boldsymbol{\mathcal{B}}^{(2)}}^{\mathbf{f}^{(2)}}\rrbracket}\left(
\mathrm{pr}^{\llbracket\cdot\rrbracket\mathsf{p}}_{\boldsymbol{\mathcal{B}}^{(2)}, \varphi(s)} \left(
x
\right)
\right).
$$

The proof of this case is identical to that of the $0$-source.

Hence, $\mathrm{pr}^{\llbracket\cdot\rrbracket\mathsf{p}}_{\boldsymbol{\mathcal{B}}^{(2)}, \varphi}$ is compatible with the $0$-target operation.

{\sffamily The mapping $\mathrm{pr}^{\llbracket\cdot\rrbracket\mathsf{p}}_{\boldsymbol{\mathcal{B}}^{(2)}, \varphi}$ is compatible with the $0$-composition.}

Let $s$ be a sort in $S$ and let us consider the $0$-composition operation symbol $\circ_{s}^{0}$ in $\Sigma_{ss,s}^{\boldsymbol{\mathcal{A}}^{(2)}}$. Let $x$ and $y$ be two elements in $\T_{\boldsymbol{\mathcal{E}}^{\boldsymbol{\mathcal{B}}^{(2)}}}(\mathbf{Pth}_{\boldsymbol{\mathcal{B}}^{(2)}})_{\varphi(s)}$ such that
$$
\mathrm{sc}_{\varphi(s)}^{(0,2)\mathbf{T}_{\boldsymbol{\mathcal{E}}^{\boldsymbol{\mathcal{B}}^{(2)}}}^{\mathbf{f}^{(2)}}(\mathbf{Pth}_{\boldsymbol{\mathcal{B}}^{(2)}})}\left(
y
\right)
=
\mathrm{tg}_{\varphi(s)}^{(0,2)\mathbf{T}_{\boldsymbol{\mathcal{E}}^{\boldsymbol{\mathcal{B}}^{(2)}}}^{\mathbf{f}^{(2)}}(\mathbf{Pth}_{\boldsymbol{\mathcal{B}}^{(2)}})}\left(
x
\right).
$$
Then the following equality holds
$$
\mathrm{pr}^{\llbracket\cdot\rrbracket\mathsf{p}}_{\boldsymbol{\mathcal{B}}^{(2)}, \varphi(s)} \left(
y
\circ_{s}^{0\mathbf{T}_{\boldsymbol{\mathcal{E}}^{\boldsymbol{\mathcal{B}}^{(2)}}}^{\mathbf{f}^{(2)}}(\mathbf{Pth}_{\boldsymbol{\mathcal{B}}^{(2)}})}\left(
x
\right)
\right)
=
\mathrm{pr}^{\llbracket\cdot\rrbracket\mathsf{p}}_{\boldsymbol{\mathcal{B}}^{(2)}, \varphi(s)} \left(
y
\right)
\circ_{s}^{0\llbracket\mathbf{Pth}_{\boldsymbol{\mathcal{B}}^{(2)}}^{\mathbf{f}^{(2)}}\rrbracket}
\mathrm{pr}^{\llbracket\cdot\rrbracket\mathsf{p}}_{\boldsymbol{\mathcal{B}}^{(2)}, \varphi(s)} \left(
x
\right).
$$

The proof of this case is identical to that of the $0$-source.

Hence, $\mathrm{pr}^{\llbracket\cdot\rrbracket\mathsf{p}}_{\boldsymbol{\mathcal{B}}^{(2)}, \varphi}$ is compatible with the $0$-composition operation.

{\sffamily The mapping $\mathrm{pr}^{\llbracket\cdot\rrbracket\mathsf{p}}_{\boldsymbol{\mathcal{B}}^{(2)}, \varphi}$ is compatible with the second-order rewrite rules.}

Let $s$ be a sort in $S$ and $\mathfrak{p}^{(2)}$ a rewrite rule in $\mathcal{A}_{s}^{(2)}$. Thus, the following chain of equalities holds
\begin{flushleft}
$
\mathrm{pr}^{\llbracket\cdot\rrbracket\mathsf{p}}_{\boldsymbol{\mathcal{B}}^{(2)}, \varphi(s)} \left(
\mathfrak{p}^{(2)\mathbf{T}_{\boldsymbol{\mathcal{E}}^{\boldsymbol{\mathcal{B}}^{(2)}}}^{\mathbf{f}^{(2)}}(\mathbf{Pth}_{\boldsymbol{\mathcal{B}}^{(2)}})}
\right)
$
\allowdisplaybreaks
\begin{align*}
&=
\mathrm{pr}^{\llbracket\cdot\rrbracket\mathsf{p}}_{\boldsymbol{\mathcal{B}}^{(2)}, \varphi(s)} \left(
\mathrm{pr}^{\equiv^{\llbracket2\rrbracket}}_{\boldsymbol{\mathcal{B}}^{(2)}, \varphi(s)} \circ \mathrm{ip}^{(2,Y)@}_{\boldsymbol{\mathcal{B}}^{(2)}, \varphi(s)} \circ \mathrm{CH}^{(2)}_{\boldsymbol{\mathcal{B}}^{(2)}, \varphi(s)} \left(
f^{(2)}_{s}\left(
\mathfrak{p}^{(2)}
\right)
\right)
\right)
\tag{1}
\\
&=
\mathrm{pr}^{\llbracket\cdot\rrbracket\mathsf{p}}_{\boldsymbol{\mathcal{B}}^{(2)}, \varphi(s)} \left(
(\mathrm{pr}^{\equiv^{\llbracket2\rrbracket}}_{\boldsymbol{\mathcal{B}}^{(2)}} \circ \mathrm{ip}^{(2,Y)@}_{\boldsymbol{\mathcal{B}}^{(2)}} \circ \mathrm{CH}^{(2)}_{\boldsymbol{\mathcal{B}}^{(2)}})_{\varphi(s)} \left(
f^{(2)}_{s}\left(
\mathfrak{p}^{(2)}
\right)
\right)
\right)
\tag{2}
\\
&=
\resizebox{.89\textwidth}{!}{%
$
\mathrm{pr}^{\llbracket\cdot\rrbracket\mathsf{p}}_{\boldsymbol{\mathcal{B}}^{(2)}, \varphi(s)} \left(
(\mathrm{pr}^{\equiv^{\llbracket2\rrbracket}}_{\boldsymbol{\mathcal{B}}^{(2)}} \circ \mathrm{ip}^{(2,Y)@}_{\boldsymbol{\mathcal{B}}^{(2)}} \circ \mathrm{CH}^{(2)}_{\boldsymbol{\mathcal{B}}^{(2)}})_{\varphi(s)}^{\natural} \left(
\left\llbracket
f^{(2)}_{s}\left(
\mathfrak{p}^{(2)}
\right)
\right\rrbracket_{\varphi(s)}
\right)
\right)
$
}
\tag{3}
\\
&=
\left\llbracket
f^{(2)}_{s}\left(
\mathfrak{p}^{(2)}
\right)
\right\rrbracket_{\varphi(s)}
\tag{4}
\\
&=
\mathfrak{p}^{(2)\llbracket\mathbf{Pth}_{\boldsymbol{\mathcal{B}}^{(2)}}^{\mathbf{f}^{(2)}}\rrbracket}.
\end{align*}
\end{flushleft}

The first equality unravels the interpretation of the constant $\mathfrak{p}^{(2)}$ in the partial $\Sigma^{\boldsymbol{\mathcal{A}}^{(2)}}$-algebra $\mathbf{T}_{\boldsymbol{\mathcal{E}}^{\boldsymbol{\mathcal{B}}^{(2)}}}^{\mathbf{f}^{(2)}}(\mathbf{Pth}_{\boldsymbol{\mathcal{B}}^{(2)}})$, introduced in Proposition~\ref{PDFreeBDCatAlg};
the second equality recovers the definition of the $\varphi(s)$ component of a composition of many-sorted mappings;
the third equality follows from Remark~\ref{RDQPUniv};
The fourth equality follows from Theorem~\ref{TDPthFree};
finally, the last equality recovers the interpretation of the constant $\mathfrak{p}^{(2)}$ in the partial $\Sigma^{\boldsymbol{\mathcal{A}}^{(2)}}$-algebra $\llbracket\mathbf{Pth}_{\boldsymbol{\mathcal{B}}^{(2)}}^{\mathbf{f}^{(2)}}\rrbracket$, introduced in Proposition~\ref{PDQPthBDCatAlg}.

Hence, $\mathrm{pr}^{\llbracket\cdot\rrbracket\mathsf{p}}_{\boldsymbol{\mathcal{B}}^{(2)}, \varphi}$ is compatible with the second-order rewrite rules.

{\sffamily The mapping $\mathrm{pr}^{\llbracket\cdot\rrbracket\mathsf{p}}_{\boldsymbol{\mathcal{B}}^{(2)}, \varphi}$ is compatible with the $1$-source.}

Let $s$ be a sort in $S$ and let us consider the $1$-source operation symbol $\mathrm{sc}_{s}^{1}$ in $\Sigma^{\boldsymbol{\mathcal{A}}^{(2)}}_{s,s}$. Let $x$ be an element in $\T_{\boldsymbol{\mathcal{E}}^{\boldsymbol{\mathcal{B}}^{(2)}}}(\mathbf{Pth}_{\boldsymbol{\mathcal{B}}^{(2)}})_{\varphi(s)}$.

The following chain of equalities holds
\allowdisplaybreaks
\begin{align*}
\mathrm{pr}^{\llbracket\cdot\rrbracket\mathsf{p}}_{\boldsymbol{\mathcal{B}}^{(2)}, \varphi(s)} \left(
\mathrm{sc}_{s}^{1\mathbf{T}_{\boldsymbol{\mathcal{E}}^{\boldsymbol{\mathcal{B}}^{(2)}}}^{\mathbf{f}^{(2)}}(\mathbf{Pth}_{\boldsymbol{\mathcal{B}}^{(2)}})}\left(
x
\right)
\right)
&=
\mathrm{pr}^{\llbracket\cdot\rrbracket\mathsf{p}}_{\boldsymbol{\mathcal{B}}^{(2)}, \varphi(s)} \left(
\mathrm{sc}_{\varphi(s)}^{1\mathbf{T}_{\boldsymbol{\mathcal{E}}^{\boldsymbol{\mathcal{B}}^{(2)}}}(\mathbf{Pth}_{\boldsymbol{\mathcal{B}}^{(2)}})}\left(
x
\right)
\right)
\tag{1}
\\
&=
\mathrm{sc}_{\varphi(s)}^{1\llbracket\mathbf{Pth}_{\boldsymbol{\mathcal{B}}^{(2)}}\rrbracket}\left(
\mathrm{pr}^{\llbracket\cdot\rrbracket\mathsf{p}}_{\boldsymbol{\mathcal{B}}^{(2)}, \varphi(s)} \left(
x
\right)
\right)
\tag{2}
\\
&=
\mathrm{sc}_{s}^{1\llbracket\mathbf{Pth}_{\boldsymbol{\mathcal{B}}^{(2)}}^{\mathbf{f}^{(2)}}\rrbracket}\left(
\mathrm{pr}^{\llbracket\cdot\rrbracket\mathsf{p}}_{\boldsymbol{\mathcal{B}}^{(2)}, \varphi(s)} \left(
x
\right)
\right).
\tag{3}
\end{align*}
The first equality unravels the interpretation of the operation symbol $\mathrm{sc}_{s}^{1}$ in the partial $\Sigma^{\boldsymbol{\mathcal{A}}^{(2)}}$-algebra $\mathbf{T}_{\boldsymbol{\mathcal{E}}^{\boldsymbol{\mathcal{B}}^{(2)}}}^{\mathbf{f}^{(2)}}(\mathbf{Pth}_{\boldsymbol{\mathcal{B}}^{(2)}})$, introduced in Proposition~\ref{PDFreeBDCatAlg};
the second equality follows from the fact that, according to Proposition~\ref{CDVarPr}, $\mathrm{pr}^{\llbracket\cdot\rrbracket\mathsf{p}}_{\boldsymbol{\mathcal{B}}^{(2)}}$ is a $\Lambda^{\boldsymbol{\mathcal{B}}^{(2)}}$-homomorphism from $\mathbf{T}_{\boldsymbol{\mathcal{E}}^{\boldsymbol{\mathcal{B}}^{(2)}}}(\mathbf{Pth}_{\boldsymbol{\mathcal{B}}^{(2)}})$ to $\llbracket\mathbf{Pth}_{\boldsymbol{\mathcal{B}}^{(2)}}\rrbracket$;
finally, the last equality recovers the interpretation of the operation symbol $\mathrm{sc}_{s}^{1}$ in the partial $\Sigma^{\boldsymbol{\mathcal{A}}^{(2)}}$-algebra $\llbracket\mathbf{Pth}_{\boldsymbol{\mathcal{B}}^{(2)}}^{\mathbf{f}^{(2)}}\rrbracket$, introduced in Proposition~\ref{PDQPthBDCatAlg}.
Hence, $\mathrm{pr}^{\llbracket\cdot\rrbracket\mathsf{p}}_{\boldsymbol{\mathcal{B}}^{(2)}, \varphi}$ is compatible with the $1$-source operation.

{\sffamily The mapping $\mathrm{pr}^{\llbracket\cdot\rrbracket\mathsf{p}}_{\boldsymbol{\mathcal{B}}^{(2)}, \varphi}$ is compatible with the $1$-target.}

Let $s$ be a sort in $S$ and let us consider the $1$-target operation symbol $\mathrm{tg}_{s}^{1}$ in $\Sigma^{\boldsymbol{\mathcal{A}}^{(2)}}_{s,s}$. Let $x$ be an element in $\T_{\boldsymbol{\mathcal{E}}^{\boldsymbol{\mathcal{B}}^{(2)}}}(\mathbf{Pth}_{\boldsymbol{\mathcal{B}}^{(2)}})_{\varphi(s)}$, then the following equality holds
$$
\mathrm{pr}^{\llbracket\cdot\rrbracket\mathsf{p}}_{\boldsymbol{\mathcal{B}}^{(2)}, \varphi(s)} \left(
\mathrm{tg}_{s}^{1\mathbf{T}_{\boldsymbol{\mathcal{E}}^{\boldsymbol{\mathcal{B}}^{(2)}}}^{\mathbf{f}^{(2)}}(\mathbf{Pth}_{\boldsymbol{\mathcal{B}}^{(2)}})}\left(
x
\right)
\right)
=
\mathrm{tg}_{s}^{1\llbracket\mathbf{Pth}_{\boldsymbol{\mathcal{B}}^{(2)}}^{\mathbf{f}^{(2)}}\rrbracket}\left(
\mathrm{pr}^{\llbracket\cdot\rrbracket\mathsf{p}}_{\boldsymbol{\mathcal{B}}^{(2)}, \varphi(s)} \left(
x
\right)
\right).
$$

The proof of this case is identical to that of the $1$-source.

Hence, $\mathrm{pr}^{\llbracket\cdot\rrbracket\mathsf{p}}_{\boldsymbol{\mathcal{B}}^{(2)}, \varphi}$ is compatible with the $1$-target operation.

{\sffamily The mapping $\mathrm{pr}^{\llbracket\cdot\rrbracket\mathsf{p}}_{\boldsymbol{\mathcal{B}}^{(2)}, \varphi}$ is compatible with the $1$-composition.}

Let $s$ be a sort in $S$ and let us consider the $1$-composition operation symbol $\circ_{s}^{1}$ in $\Sigma_{ss,s}^{\boldsymbol{\mathcal{A}}^{(2)}}$. Let $x$ and $y$ be two elements in $\T_{\boldsymbol{\mathcal{E}}^{\boldsymbol{\mathcal{B}}^{(2)}}}(\mathbf{Pth}_{\boldsymbol{\mathcal{B}}^{(2)}})_{\varphi(s)}$ such that
$$
\mathrm{sc}_{\varphi(s)}^{(1,2)\mathbf{T}_{\boldsymbol{\mathcal{E}}^{\boldsymbol{\mathcal{B}}^{(2)}}}^{\mathbf{f}^{(2)}}(\mathbf{Pth}_{\boldsymbol{\mathcal{B}}^{(2)}})}\left(
y
\right)
=
\mathrm{tg}_{\varphi(s)}^{(1,2)\mathbf{T}_{\boldsymbol{\mathcal{E}}^{\boldsymbol{\mathcal{B}}^{(2)}}}^{\mathbf{f}^{(2)}}(\mathbf{Pth}_{\boldsymbol{\mathcal{B}}^{(2)}})}\left(
x
\right).
$$
Then the following equality holds
$$
\mathrm{pr}^{\llbracket\cdot\rrbracket\mathsf{p}}_{\boldsymbol{\mathcal{B}}^{(2)}, \varphi(s)} \left(
y
\circ_{s}^{1\mathbf{T}_{\boldsymbol{\mathcal{E}}^{\boldsymbol{\mathcal{B}}^{(2)}}}^{\mathbf{f}^{(2)}}(\mathbf{Pth}_{\boldsymbol{\mathcal{B}}^{(2)}})}\left(
x
\right)
\right)
=
\mathrm{pr}^{\llbracket\cdot\rrbracket\mathsf{p}}_{\boldsymbol{\mathcal{B}}^{(2)}, \varphi(s)} \left(
y
\right)
\circ_{s}^{1\llbracket\mathbf{Pth}_{\boldsymbol{\mathcal{B}}^{(2)}}^{\mathbf{f}^{(2)}}\rrbracket}
\mathrm{pr}^{\llbracket\cdot\rrbracket\mathsf{p}}_{\boldsymbol{\mathcal{B}}^{(2)}, \varphi(s)} \left(
x
\right).
$$

The proof of this case is identical to that of the $1$-source.

Hence, $\mathrm{pr}^{\llbracket\cdot\rrbracket\mathsf{p}}_{\boldsymbol{\mathcal{B}}^{(2)}, \varphi}$ is compatible with the $1$-composition operation.

This completes the proof.
\end{proof}

\begin{theorem}\label{TDPthBFreeB}
The partial $\Sigma^{\boldsymbol{\mathcal{A}}^{(2)}}$-algebras $\mathbf{T}_{\boldsymbol{\mathcal{E}}^{\boldsymbol{\mathcal{B}}^{(2)}}}^{\mathbf{f}^{(2)}}(\mathbf{Pth}_{\boldsymbol{\mathcal{B}}^{(2)}})$ and $\llbracket\mathbf{Pth}_{\boldsymbol{\mathcal{B}}^{(2)}}^{\mathbf{f}^{(2)}}\rrbracket$ are isomorphic.
\end{theorem}

\begin{proof}
Note that it follows from Theorem~\ref{TDPthFree} that $(\mathrm{pr}^{\equiv^{\llbracket2\rrbracket}}_{\boldsymbol{\mathcal{B}}^{(2)}} \circ \mathrm{ip}^{(2,Y)@}_{\boldsymbol{\mathcal{B}}^{(2)}} \circ \mathrm{CH}^{(2)}_{\boldsymbol{\mathcal{B}}^{(2)}})^{\natural}$ and $\mathrm{pr}^{\llbracket\cdot\rrbracket\mathsf{p}}_{\boldsymbol{\mathcal{B}}^{(2)}}$ are a pair of inverse $\Lambda^{\boldsymbol{\mathcal{B}}^{(2)}}$-isomorphisms.

Thus, for every sort $s$ in $S$, the following chain of equalities holds
\begin{flushleft}
$
(\mathrm{pr}^{\equiv^{\llbracket2\rrbracket}}_{\boldsymbol{\mathcal{B}}^{(2)}} \circ \mathrm{ip}^{(2,Y)@}_{\boldsymbol{\mathcal{B}}^{(2)}} \circ \mathrm{CH}^{(2)}_{\boldsymbol{\mathcal{B}}^{(2)}})^{\natural}_{\varphi(s)}
\circ
\mathrm{pr}^{\llbracket\cdot\rrbracket\mathsf{p}}_{\boldsymbol{\mathcal{B}}^{(2)}, \varphi(s)}
$
\begin{align*}
&=
\left(
(\mathrm{pr}^{\equiv^{\llbracket2\rrbracket}}_{\boldsymbol{\mathcal{B}}^{(2)}} \circ \mathrm{ip}^{(2,Y)@}_{\boldsymbol{\mathcal{B}}^{(2)}} \circ \mathrm{CH}^{(2)}_{\boldsymbol{\mathcal{B}}^{(2)}})^{\natural}
\circ
\mathrm{pr}^{\llbracket\cdot\rrbracket\mathsf{p}}_{\boldsymbol{\mathcal{B}}^{(2)}}
\right)_{\varphi(s)}
\tag{1}
\\
&=
\left(
\T_{\boldsymbol{\mathcal{E}}^{\boldsymbol{\mathcal{B}}^{(2)}}}(\mathbf{Pth}_{\boldsymbol{\mathcal{B}}^{(2)}})
\right)_{\varphi(s)}
\tag{2}
\end{align*}

Likewise, for every sort $s$ in $S$, the following chain of equalities holds
\begin{flushleft}
$
\mathrm{pr}^{\llbracket\cdot\rrbracket\mathsf{p}}_{\boldsymbol{\mathcal{B}}^{(2)}, \varphi(s)}
\circ
(\mathrm{pr}^{\equiv^{\llbracket2\rrbracket}}_{\boldsymbol{\mathcal{B}}^{(2)}} \circ \mathrm{ip}^{(2,Y)@}_{\boldsymbol{\mathcal{B}}^{(2)}} \circ \mathrm{CH}^{(2)}_{\boldsymbol{\mathcal{B}}^{(2)}})^{\natural}_{\varphi(s)}
$
\end{flushleft}
\begin{align*}
&=
\left(
\mathrm{pr}^{\llbracket\cdot\rrbracket\mathsf{p}}_{\boldsymbol{\mathcal{B}}^{(2)}}
\circ
(\mathrm{pr}^{\equiv^{\llbracket2\rrbracket}}_{\boldsymbol{\mathcal{B}}^{(2)}} \circ \mathrm{ip}^{(2,Y)@}_{\boldsymbol{\mathcal{B}}^{(2)}} \circ \mathrm{CH}^{(2)}_{\boldsymbol{\mathcal{B}}^{(2)}})^{\natural}
\right)_{\varphi(s)}
\tag{1}
\\
&=
\left(
\llbracket\mathrm{Pth}_{\boldsymbol{\mathcal{B}}^{(2)}}\rrbracket
\right)_{\varphi(s)}
\tag{2}
\end{align*}
\end{flushleft}
\end{proof}

\begin{corollary}\label{TDPTBFreeB}
The partial $\Sigma^{\boldsymbol{\mathcal{A}}^{(2)}}$-algebras $\mathbf{T}_{\boldsymbol{\mathcal{E}}^{\boldsymbol{\mathcal{B}}^{(2)}}}^{\mathbf{f}^{(2)}}(\mathbf{Pth}_{\boldsymbol{\mathcal{B}}^{(2)}})$ and $\llbracket\mathbf{PT}_{\boldsymbol{\mathcal{B}}^{(2)}}^{\mathbf{f}^{(2)}}\rrbracket$ are isomorphic.
\end{corollary}			
\chapter{Second-order quotient path extension mapping}\label{S3G}

The aim of this chapter is to define, given a second-order morphism $\mathbf{f}^{(2)}$ and its second-order path extension $f^{(2)\flat}$, a mapping from the $S$-sorted set $\llbracket \mathrm{Pth}_{\boldsymbol{\mathcal{A}}^{(2)}} \rrbracket$ to the $S$-sorted set $\llbracket \mathrm{Pth}_{\boldsymbol{\mathcal{B}}^{(2)}} \rrbracket_{\varphi}$ that we will call the second-order quotient path-extension mapping denoted by $f^{\llbracket 2 \rrbracket @}$ which, by construction, is a $\Sigma^{\boldsymbol{\mathcal{A}}^{(2)}}$-homomorphism from $\mathbf{Pth}_{\boldsymbol{\mathcal{A}}^{(2)}}$ to $\mathbf{Pth}_{\boldsymbol{\mathcal{B}}^{(2)}}^{\mathbf{f}^{(2)}}$.

In order to construct the desired mapping, we begin by proving that the partial $\Sigma^{\boldsymbol{\mathcal{A}}^{(2)}}$-algebra $\llbracket \mathbf{Pth}_{\boldsymbol{\mathcal{B}}^{(2)}}^{\mathbf{f}^{(2)}} \rrbracket$ satisfies axioms A8 and B8 defining the QE-variety $\mathcal{V}(\boldsymbol{\mathcal{E}}^{\boldsymbol{\mathcal{A}}^{(2)}})$, introduced in Definition~\ref{DDVar}.

\begin{proposition}
\label{PDQPthBVarA8}
Let $(\mathbf{s}, s)$ an element of $S^{\star} \times S$, $\sigma$ an operation symbol in $\Sigma_{\mathbf{s}, s}$ and $( \llbracket\mathfrak{P}_{j}^{(2)}\rrbracket_{\varphi(s_{j})} )_{j\in\bb{\mathbf{s}}}$ and $( \llbracket\mathfrak{Q}_{j}^{(2)}\rrbracket_{\varphi(s_{j})} )_{j\in\bb{\mathbf{s}}}$ be two family of path classes in $\llbracket\mathrm{Pth}_{\boldsymbol{\mathcal{B}}^{(2)}}\rrbracket_{\varphi^{\star}(\mathbf{s})}$ such that, for every $j \in \bb{\mathbf{s}}$, 
$$
\mathrm{sc}_{s_{j}}^{0\llbracket\mathbf{Pth}_{\boldsymbol{\mathcal{B}}^{(2)}}^{\mathbf{f}^{(2)}}\rrbracket} \left(
\left\llbracket
\mathfrak{Q}_{j}^{(2)}
\right\rrbracket_{s_{j}}
\right)
=
\mathrm{tg}_{s_{j}}^{0\llbracket\mathbf{Pth}_{\boldsymbol{\mathcal{B}}^{(1)}}^{\mathbf{f}^{(1)}}\rrbracket} \left(
\left\llbracket
\mathfrak{P}_{j}^{(2)}
\right\rrbracket_{s_{j}}
\right).
$$
Then the following equality holds
\begin{multline*}
\sigma^{\llbracket\mathbf{Pth}_{\boldsymbol{\mathcal{B}}^{(1)}}^{\mathbf{f}^{(1)}}\rrbracket} \left(
\left(
\left\llbracket
\mathfrak{Q}_{j}^{(2)}
\right\rrbracket_{s_{j}}
\circ_{s_{j}}^{0\llbracket\mathbf{Pth}_{\boldsymbol{\mathcal{B}}^{(1)}}^{\mathbf{f}^{(1)}}\rrbracket}
\left\llbracket
\mathfrak{P}_{j}^{(2)}
\right\rrbracket_{s_{j}}
\right)_{j \in \bb{\mathbf{s}}}
\right)
\\
=
\sigma^{\llbracket\mathbf{Pth}_{\boldsymbol{\mathcal{B}}^{(2)}}^{\mathbf{f}^{(2)}}\rrbracket} \left(
\left(
\left\llbracket
\mathfrak{Q}_{j}^{(2)}
\right\rrbracket_{s_{j}}
\right)_{j \in \bb{\mathbf{s}}}
\right)
\circ_{s}^{0\llbracket\mathbf{Pth}_{\boldsymbol{\mathcal{B}}^{(1)}}^{\mathbf{f}^{(1)}}\rrbracket}
\sigma^{\llbracket\mathbf{Pth}_{\boldsymbol{\mathcal{B}}^{(1)}}^{\mathbf{f}^{(1)}}\rrbracket} \left(
\left(
\left\llbracket
\mathfrak{P}_{j}^{(2)}
\right\rrbracket_{s_{j}}
\right)_{j \in \bb{\mathbf{s}}}
\right).
\end{multline*}
\end{proposition}

\begin{proof}
Let $(\mathbf{s}, s)$ an element of $S^{\star} \times S$, $\sigma$ an operation symbol in $\Sigma_{\mathbf{s}, s}$ and $( \llbracket\mathfrak{P}^{(2)}_{j}\rrbracket_{\varphi(s_{j})} )_{j\in\bb{\mathbf{s}}}$ and $( \llbracket\mathfrak{Q}^{(2)}_{j}\rrbracket_{\varphi(s_{j})} )_{j\in\bb{\mathbf{s}}}$ be two family of path classes in $\llbracket\mathrm{Pth}_{\boldsymbol{\mathcal{B}}^{(2)}}\rrbracket_{\varphi^{\star}(\mathbf{s})}$ such that, for every $j \in \bb{\mathbf{s}}$, 
$$
\mathrm{sc}_{s_{j}}^{0\llbracket\mathbf{Pth}_{\boldsymbol{\mathcal{B}}^{(2)}}^{\mathbf{f}^{(2)}}\rrbracket} \left(
\left\llbracket
\mathfrak{Q}_{j}^{(2)}
\right\rrbracket_{s_{j}}
\right)
=
\mathrm{tg}_{s_{j}}^{0\llbracket\mathbf{Pth}_{\boldsymbol{\mathcal{B}}^{(1)}}^{\mathbf{f}^{(1)}}\rrbracket} \left(
\left\llbracket
\mathfrak{P}_{j}^{(2)}
\right\rrbracket_{s_{j}}
\right).
$$
Unraveling the definition of the interpretation of the operation symbols $\sigma$ and $\circ_{s}^{0}$ in the $\Sigma^{\boldsymbol{\mathcal{A}}^{(2)}}$-algebras $\llbracket\mathbf{Pth}_{\boldsymbol{\mathcal{B}}^{(2)}}^{\mathbf{f}^{(2)}}\rrbracket$ and $\mathbf{Pth}_{\boldsymbol{\mathcal{B}}^{(2)}}^{\mathbf{f}^{(2)}}$ introduced in Propositions~\ref{PDQPthBDCatAlg} and \ref{PDPthBCatAlg}, respectively, the desired equality is equivalent to
\begin{multline*}
\left\llbracket
\sigma^{\mathbf{c}_{\mathfrak{d}}^{\ast}(\mathbf{Pth}_{\boldsymbol{\mathcal{B}}^{(2)}}^{(0,2)})} \left(
\left(
\mathfrak{Q}_{j}^{(2)}
\circ_{\varphi(s_{j})}^{0\mathbf{Pth}_{\boldsymbol{\mathcal{B}}^{(2)}}}
\mathfrak{P}_{j}^{(2)}
\right)_{j\in\bb{\mathbf{s}}}
\right)
\right\rrbracket_{\varphi(s)}
\\
=
\left\llbracket
\sigma^{\mathbf{c}_{\mathfrak{d}}^{\ast}(\mathbf{Pth}_{\boldsymbol{\mathcal{B}}^{(2)}}^{(0,2)})} \left(
\left(
\mathfrak{Q}_{j}^{(2)}
\right)_{j \in \bb{\mathbf{s}}}
\right)
\circ_{\varphi(s)}^{0\mathbf{Pth}_{\boldsymbol{\mathcal{B}}^{(2)}}}
\sigma^{\mathbf{c}_{\mathfrak{d}}^{\ast}(\mathbf{Pth}_{\boldsymbol{\mathcal{B}}^{(2)}}^{(0,2)})} \left(
\left(
\mathfrak{P}_{j}^{(2)}
\right)_{j \in \bb{\mathbf{s}}}
\right)
\right\rrbracket_{\varphi(s)}.
\tag{\dag}
\end{multline*}

Let us recall that, the interpretation of the operation symbol $\sigma$ in the $\Sigma$-algebra $\mathbf{c}_{\mathfrak{d}}^{\ast}(\mathbf{Pth}_{\boldsymbol{\mathcal{B}}^{(2)}}^{(0,2)})$ is given by the derived operation on $\mathbf{Pth}_{\boldsymbol{\mathcal{B}}^{(2)}}^{(0,2)}$ determined by $c_{\mathbf{s}, s}(\sigma)$ which following Remark~\ref{Rdop}, we denote by $c(\sigma)^{\mathbf{Pth}_{\boldsymbol{\mathcal{B}}^{(2)}}^{(0,2)}}$. For the details about derived operations, see Definition~\ref{wdop}.

Thus, we prove that, for every $(\mathbf{t}, t)$ in $T^{\star} \times T$, every $P$ in $T_{\Lambda}(\vs\mathbf{t})_{t}$ and every pair of families $(\mathfrak{P}_{j}^{(2)})_{j \in \bb{\mathbf{t}}}$ and $(\mathfrak{Q}_{j}^{(2)})_{j \in \bb{\mathbf{t}}}$ in $\mathrm{Pth}_{\boldsymbol{\mathcal{B}}^{(2)}, \mathbf{t}}$,
\begin{multline*}
\left\llbracket
P^{\mathbf{Pth}_{\boldsymbol{\mathcal{B}}^{(2)}}^{(0,2)}} \left(
\left(
\mathfrak{Q}_{j}^{(2)}
\circ_{t_{j}}^{0\mathbf{Pth}_{\boldsymbol{\mathcal{B}}^{(2)}}}
\mathfrak{P}_{j}^{(2)}
\right)_{j\in\bb{\mathbf{t}}}
\right)
\right\rrbracket_{t}
\\
=
\left\llbracket
P^{\mathbf{Pth}_{\boldsymbol{\mathcal{B}}^{(2)}}^{(0,2)}} \left(
\left(
\mathfrak{Q}_{j}^{(2)}
\right)_{j \in \bb{\mathbf{t}}}
\right)
\circ_{t}^{0\mathbf{Pth}_{\boldsymbol{\mathcal{B}}^{(2)}}}
P^{\mathbf{Pth}_{\boldsymbol{\mathcal{B}}^{(2)}}^{(0,2)}} \left(
\left(
\mathfrak{P}_{j}^{(2)}
\right)_{j \in \bb{\mathbf{t}}}
\right)
\right\rrbracket_{t}.
\end{multline*}
Therefore, Equation~(\dag) follows as a consequence.

By algebraic induction on the complexity of $P$.

{\sffamily Base step of the induction}

If $P$ is a variable $v_{i}^{t}$, thus, $i \in \bb{t}$ and $t_{i} = t$, then following Remark~\ref{Rdop}
$$
(v_{i}^{t})^{\mathbf{Pth}_{\boldsymbol{\mathcal{B}}^{(2)}}^{(0,2)}}
=
\mathrm{d}_{\mathbf{t}, t}^{\mathbf{Pth}_{\boldsymbol{\mathcal{B}}^{(2)}}^{(0,2)}} (v_{i}^{t})
=
\mathrm{d}_{\mathbf{t}, t}^{\mathbf{Pth}_{\boldsymbol{\mathcal{B}}^{(2)}}^{(0,2)}} ( \eta^{\vs\mathbf{t}} (v_{i}^{t})) 
=
\mathrm{p}_{\mathbf{t}, t}^{\mathbf{Pth}_{\boldsymbol{\mathcal{B}}^{(2)}}^{(0,2)}} (v_{i}^{t})
=
\mathrm{pr}_{\mathbf{t}, i}^{\mathrm{Pth}_{\boldsymbol{\mathcal{B}}^{(2)}}}.
$$

The following chain of equalities holds
\allowdisplaybreaks
\begin{flushleft}
$
P^{\mathbf{Pth}_{\boldsymbol{\mathcal{B}}^{(2)}}^{(0,2)}} \left(
\left(
\mathfrak{Q}_{j}^{(2)}
\circ_{t_{j}}^{0\mathbf{Pth}_{\boldsymbol{\mathcal{B}}^{(2)}}}
\mathfrak{P}_{j}^{(2)}
\right)_{j\in\bb{\mathbf{t}}}
\right)
$
\begin{align*}
&=
(v_{i}^{t})^{\mathbf{Pth}_{\boldsymbol{\mathcal{B}}^{(2)}}^{(0,2)}} \left(
\left(
\mathfrak{Q}_{j}^{(2)}
\circ_{t_{j}}^{0\mathbf{Pth}_{\boldsymbol{\mathcal{B}}^{(2)}}}
\mathfrak{P}_{j}^{(2)}
\right)_{j\in\bb{\mathbf{t}}}
\right)
\tag{1}
\\
&=
\mathrm{pr}_{\mathbf{t}, i}^{\mathrm{Pth}_{\boldsymbol{\mathcal{B}}^{(2)}}} \left(
\left(
\mathfrak{Q}_{j}^{(2)}
\circ_{t_{j}}^{0\mathbf{Pth}_{\boldsymbol{\mathcal{B}}^{(2)}}}
\mathfrak{P}_{j}^{(2)}
\right)_{j\in\bb{\mathbf{t}}}
\right)
\tag{2}
\\
&=
\mathfrak{Q}_{i}^{(2)}
\circ_{t_{i}}^{0\mathbf{Pth}_{\boldsymbol{\mathcal{B}}^{(2)}}}
\mathfrak{P}_{i}^{(2)}
\tag{3}
\\
&=
\mathfrak{Q}_{i}^{(2)}
\circ_{t}^{0\mathbf{Pth}_{\boldsymbol{\mathcal{B}}^{(2)}}}
\mathfrak{P}_{i}^{(2)}
\tag{4}
\\
&=
\mathrm{pr}_{\mathbf{t}, i}^{\mathrm{Pth}_{\boldsymbol{\mathcal{B}}^{(2)}}} \left(
\left(
\mathfrak{Q}_{j}^{(2)}
\right)_{j \in \bb{\mathbf{t}}}
\right)
\circ_{t}^{0\mathbf{Pth}_{\boldsymbol{\mathcal{B}}^{(2)}}}
\mathrm{pr}_{\mathbf{t}, i}^{\mathrm{Pth}_{\boldsymbol{\mathcal{B}}^{(2)}}} \left(
\left(
\mathfrak{P}_{j}^{(2)}
\right)_{j\in\bb{\mathbf{t}}}
\right)
\tag{5}
\\
&=
(v_{i}^{t})^{\mathbf{Pth}_{\boldsymbol{\mathcal{B}}^{(2)}}^{(0,2)}} \left(
\left(
\mathfrak{Q}_{j}^{(2)}
\right)_{j \in \bb{\mathbf{t}}}
\right)
\circ_{t}^{0\mathbf{Pth}_{\boldsymbol{\mathcal{B}}^{(2)}}}
(v_{i}^{t})^{\mathbf{Pth}_{\boldsymbol{\mathcal{B}}^{(2)}}^{(0,2)}} \left(
\left(
\mathfrak{P}_{j}^{(2)}
\right)_{j\in\bb{\mathbf{t}}}
\right)
\tag{6}
\\
&=
P^{\mathbf{Pth}_{\boldsymbol{\mathcal{B}}^{(2)}}^{(0,2)}} \left(
\left(
\mathfrak{Q}_{j}^{(2)}
\right)_{j \in \bb{\mathbf{t}}}
\right)
\circ_{t}^{0\mathbf{Pth}_{\boldsymbol{\mathcal{B}}^{(2)}}}
P^{\mathbf{Pth}_{\boldsymbol{\mathcal{B}}^{(2)}}^{(0,2)}} \left(
\left(
\mathfrak{P}_{j}^{(2)}
\right)_{j\in\bb{\mathbf{t}}}
\right).
\tag{7}
\end{align*}
\end{flushleft}

The first equality unravels the definition of $P$;
the second equality follows from the fact that $ (v_{i}^{t})^{\mathbf{Pth}_{\boldsymbol{\mathcal{B}}^{(2)}}^{(0,2)}} = \mathrm{pr}_{i}^{\mathrm{Pth}_{\boldsymbol{\mathcal{B}}^{(2)}, \mathbf{t}}} $;
the third equality unravels the definition of the $\mathbf{t}$-ary, $i$-th canonical projection from $\mathrm{Pth}_{\boldsymbol{\mathcal{B}}^{(2)}, \mathbf{t}}$ to $\mathrm{Pth}_{\boldsymbol{\mathcal{B}}^{(2)}, t_{i}}$;
the fourth equality follows from the fact that $t_{i} = t$;
the fifth equality recovers the definition of the $\mathbf{t}$-ary, $i$-th canonical projection from $\mathrm{Pth}_{\boldsymbol{\mathcal{B}}^{(2)}, \mathbf{t}}$ to $\mathrm{Pth}_{\boldsymbol{\mathcal{B}}^{(2)}, t_{i}}$;
the sixth equality follows from the fact that $ (v_{i}^{t})^{\mathbf{Pth}_{\boldsymbol{\mathcal{B}}^{(2)}}^{(0,2)}} = \mathrm{pr}_{i}^{\mathrm{Pth}_{\boldsymbol{\mathcal{B}}^{(2)}, \mathbf{t}}} $;
finally, the last equality recovers the definition of $P$.

Therefore, its $\llbracket\cdot\rrbracket$-classes coincide, i.e., equality (\dag) follows.

If $P$ is a constant symbol $\tau$ in $\Lambda_{\lambda, t}$, then the following chain of equalities holds
\allowdisplaybreaks
\begin{flushleft}
$
P^{\mathbf{Pth}_{\boldsymbol{\mathcal{B}}^{(2)}}^{(0,2)}} \left(
\left(
\mathfrak{Q}_{j}^{(2)}
\circ_{t_{j}}^{0\mathbf{Pth}_{\boldsymbol{\mathcal{B}}^{(2)}}}
\mathfrak{P}_{j}^{(2)}
\right)_{j\in\bb{\mathbf{t}}}
\right)
$
\begin{align*}
&=
\tau^{\mathbf{Pth}_{\boldsymbol{\mathcal{B}}^{(2)}}^{(0,2)}}
\tag{1}
\\
&=
\tau^{\mathbf{Pth}_{\boldsymbol{\mathcal{B}}^{(2)}}^{(0,2)}}
\circ_{t}^{0\mathbf{Pth}_{\boldsymbol{\mathcal{B}}^{(2)}}}
\tau^{\mathbf{Pth}_{\boldsymbol{\mathcal{B}}^{(2)}}^{(0,2)}}
\tag{2}
\\
&=
P^{\mathbf{Pth}_{\boldsymbol{\mathcal{B}}^{(2)}}^{(0,2)}} \left(
\left(
\mathfrak{Q}_{j}^{(2)}
\right)_{j \in \bb{\mathbf{t}}}
\right)
\circ_{t}^{0\mathbf{Pth}_{\boldsymbol{\mathcal{B}}^{(2)}}}
P^{\mathbf{Pth}_{\boldsymbol{\mathcal{B}}^{(2)}}^{(0,2)}} \left(
\left(
\mathfrak{P}_{j}^{(2)}
\right)_{j\in\bb{\mathbf{t}}}
\right).
\tag{3}
\end{align*}
\end{flushleft}

The first equality unravels the definition of $P$;
the second equality follows from the fact that, according to Proposition~\ref{PDPthCatAlg}, the interpretation of the $\tau$ operation in the many-sorted partial $\Lambda^{\boldsymbol{\mathcal{B}}^{(2)}}$-algebra $\mathbf{Pth}_{\boldsymbol{\mathcal{B}}^{(2)}}$ is given by the $(2,[1])$-identity path on $[\tau^{\mathbf{PT}_{\boldsymbol{\mathcal{B}}^{(1)}}}]$, i.e., $\tau^{\mathbf{Pth}_{\boldsymbol{\mathcal{B}}^{(2)}}} = \mathrm{ip}_{t}^{(2,[1])\sharp} \left([\tau^{\mathbf{PT}_{\boldsymbol{\mathcal{B}}^{(1)}}}]\right)$ and from the fact that, according to Proposition~\ref{PDUIp}, the $(2,[1])$-identity paths are idempotent for the $0$-composition;
finally, the last equality recovers the definition of $P$.

Therefore, its $\llbracket\cdot\rrbracket$-classes coincide, i.e., equality (\dag) follows.

{\sffamily Inductive step of the induction}

If $P$ is a term $\tau \left(\left(Q_{i}\right)_{i \in \bb{\mathbf{r}}}\right)$ for $\tau$ an operation symbol in $\lambda_{\mathbf{r}, t}$ and $(Q_{i})_{i \in \bb{\mathbf{r}}}$ a family of terms in $\T_{\Lambda}(\vs\mathbf{t})_{\mathbf{r}}$, then following Remark~\ref{Rdop}
\begin{align*}
P^{\mathbf{Pth}_{\boldsymbol{\mathcal{B}}^{(2)}}^{(0,2)}}
&=
\mathrm{d}_{\mathbf{t}, t}^{\mathbf{Pth}_{\boldsymbol{\mathcal{B}}^{(2)}}^{(0,2)}} (P)
\\
&=
\mathrm{d}_{\mathbf{t}, t}^{\mathbf{Pth}_{\boldsymbol{\mathcal{B}}^{(2)}}^{(0,2)}} \left(\tau \left(\left(Q_{i}\right)_{i \in \bb{\mathbf{r}}}\right)\right) 
\\
&=
\tau^{\mathbf{O}_{\mathbf{t}}\left(\mathbf{Pth}_{\boldsymbol{\mathcal{B}}^{(2)}}^{(0,2)}\right)} \left(
\left(
\mathrm{d}_{\mathbf{t}, r_{i}}^{\mathbf{Pth}_{\boldsymbol{\mathcal{B}}^{(2)}}^{(0,2)}}
\left(
Q_{i}
\right)
\right)_{i \in \bb{\mathbf{r}}}
\right)
\\
&=
\tau^{\mathbf{O}_{\mathbf{t}}\left(\mathbf{Pth}_{\boldsymbol{\mathcal{B}}^{(2)}}^{(0,2)}\right)} \left(
\left(
Q_{i}^{\mathbf{Pth}_{\boldsymbol{\mathcal{B}}^{(2)}}^{(0,2)}}
\right)_{i \in \bb{\mathbf{r}}}
\right).
\end{align*}

The following chain of equalities holds
\allowdisplaybreaks
\begin{flushleft}
$
\left\llbracket
P^{\mathbf{Pth}_{\boldsymbol{\mathcal{B}}^{(2)}}^{(0,2)}} \left(
\left(
\mathfrak{Q}_{j}^{(2)}
\circ_{t_{j}}^{0\mathbf{Pth}_{\boldsymbol{\mathcal{B}}^{(2)}}}
\mathfrak{P}_{j}^{(2)}
\right)_{j\in\bb{\mathbf{t}}}
\right)
\right\rrbracket_{t}
$
\begin{align*}
&=
\left\llbracket
\tau^{\mathbf{O}_{\mathbf{t}}\left(\mathbf{Pth}_{\boldsymbol{\mathcal{B}}^{(2)}}^{(0,2)}\right)} \left(
\left(
Q_{i}^{\mathbf{Pth}_{\boldsymbol{\mathcal{B}}^{(2)}}^{(0,2)}}
\right)_{i \in \bb{\mathbf{r}}}
\right)
\left(
\left(
\mathfrak{Q}_{j}^{(2)}
\circ_{t_{j}}^{0\mathbf{Pth}_{\boldsymbol{\mathcal{B}}^{(2)}}}
\mathfrak{P}_{j}^{(2)}
\right)_{j\in\bb{\mathbf{t}}}
\right)
\right\rrbracket_{t}
\tag{1}
\\
&=
\left\llbracket
\tau^{\mathbf{Pth}_{\boldsymbol{\mathcal{B}}^{(2)}}^{(0,2)}} \left(
\left(
Q_{i}^{\mathbf{Pth}_{\boldsymbol{\mathcal{B}}^{(2)}}^{(0,2)}}
\left(
\left(
\mathfrak{Q}_{j}^{(2)}
\circ_{t_{j}}^{0\mathbf{Pth}_{\boldsymbol{\mathcal{B}}^{(2)}}}
\mathfrak{P}_{j}^{(2)}
\right)_{j\in\bb{\mathbf{t}}}
\right)
\right)_{i \in \bb{\mathbf{r}}}
\right)
\right\rrbracket_{t}
\tag{2}
\\
&=
\tau^{\left\llbracket\mathbf{Pth}_{\boldsymbol{\mathcal{B}}^{(2)}}^{(0,2)}\right\rrbracket} \left(
\left(
\left\llbracket
Q_{i}^{\mathbf{Pth}_{\boldsymbol{\mathcal{B}}^{(2)}}^{(0,2)}}
\left(
\left(
\mathfrak{Q}_{j}^{(2)}
\circ_{t_{j}}^{0\mathbf{Pth}_{\boldsymbol{\mathcal{B}}^{(2)}}}
\mathfrak{P}_{j}^{(2)}
\right)_{j\in\bb{\mathbf{t}}}
\right)
\right\rrbracket_{r_{i}}
\right)_{i\in\bb{\mathbf{r}}}
\right)
\tag{3}
\\
&=
\resizebox{0.86\textwidth}{!}{%
$
\tau^{\left\llbracket\mathbf{Pth}_{\boldsymbol{\mathcal{B}}^{(2)}}^{(0,2)}\right\rrbracket} \left(
\left(
\left\llbracket
Q_{i}^{\mathbf{Pth}_{\boldsymbol{\mathcal{B}}^{(2)}}^{(0,2)}} \left(
\left(
\mathfrak{Q}_{j}^{(2)}
\right)_{j \in \bb{\mathbf{t}}}
\right)
\circ_{r_{i}}^{0\mathbf{Pth}_{\boldsymbol{\mathcal{B}}^{(2)}}}
Q_{i}^{\mathbf{Pth}_{\boldsymbol{\mathcal{B}}^{(2)}}^{(0,2)}} \left(
\left(
\mathfrak{P}_{j}^{(2)}
\right)_{j \in \bb{\mathbf{t}}}
\right)
\right\rrbracket_{r_{i}}
\right)_{i\in\bb{\mathbf{r}}}
\right)
$}
\tag{4}
\\
&=
\resizebox{0.86\textwidth}{!}{%
$
\tau^{\left\llbracket\mathbf{Pth}_{\boldsymbol{\mathcal{B}}^{(2)}}^{(0,2)}\right\rrbracket} \left(
\left(
\left\llbracket
Q_{i}^{\mathbf{Pth}_{\boldsymbol{\mathcal{B}}^{(2)}}^{(0,2)}} \left(
\left(
\mathfrak{Q}_{j}^{(2)}
\right)_{j \in \bb{\mathbf{t}}}
\right)
\right\rrbracket_{r_{i}}
\circ_{r_{i}}^{0\left\llbracket\mathbf{Pth}_{\boldsymbol{\mathcal{B}}^{(2)}}\right\rrbracket}
\left\llbracket
Q_{i}^{\mathbf{Pth}_{\boldsymbol{\mathcal{B}}^{(2)}}^{(0,2)}} \left(
\left(
\mathfrak{P}_{j}^{(2)}
\right)_{j \in \bb{\mathbf{t}}}
\right)
\right\rrbracket_{r_{i}}
\right)_{i\in\bb{\mathbf{r}}}
\right)
$}
\tag{5}
\\
&=
\tau^{\left\llbracket\mathbf{Pth}_{\boldsymbol{\mathcal{B}}^{(2)}}^{(0,2)}\right\rrbracket} \left(
\left(
\left\llbracket
Q_{i}^{\mathbf{Pth}_{\boldsymbol{\mathcal{B}}^{(2)}}^{(0,2)}} \left(
\left(
\mathfrak{Q}_{j}^{(2)}
\right)_{j \in \bb{\mathbf{t}}}
\right)
\right\rrbracket_{r_{i}}
\right)_{i\in\bb{\mathbf{r}}}
\right)
\\
&\hspace{2cm}
\circ_{t}^{0\left\llbracket\mathbf{Pth}_{\boldsymbol{\mathcal{B}}^{(2)}}\right\rrbracket}
\tau^{\left\llbracket\mathbf{Pth}_{\boldsymbol{\mathcal{B}}^{(2)}}^{(0,2)}\right\rrbracket} \left(
\left(
\left\llbracket
Q_{i}^{\mathbf{Pth}_{\boldsymbol{\mathcal{B}}^{(2)}}^{(0,2)}} \left(
\left(
\mathfrak{P}_{j}^{(2)}
\right)_{j \in \bb{\mathbf{t}}}
\right)
\right\rrbracket_{r_{i}}
\right)_{i\in\bb{\mathbf{r}}}
\right)
\tag{6}
\\
&=
\left\llbracket
\tau^{\mathbf{Pth}_{\boldsymbol{\mathcal{B}}^{(2)}}^{(0,2)}} \left(
\left(
Q_{i}^{\mathbf{Pth}_{\boldsymbol{\mathcal{B}}^{(2)}}^{(0,2)}} \left(
\left(
\mathfrak{Q}_{j}^{(2)}
\right)_{j \in \bb{\mathbf{t}}}
\right)
\right)_{i\in\bb{\mathbf{r}}}
\right)
\right\rrbracket_{t}
\\
&\hspace{2.5cm}
\circ_{t}^{0\left\llbracket\mathbf{Pth}_{\boldsymbol{\mathcal{B}}^{(2)}}\right\rrbracket}
\left\llbracket
\tau^{\mathbf{Pth}_{\boldsymbol{\mathcal{B}}^{(2)}}^{(0,2)}} \left(
\left(
Q_{i}^{\mathbf{Pth}_{\boldsymbol{\mathcal{B}}^{(2)}}^{(0,2)}} \left(
\left(
\mathfrak{P}_{j}^{(2)}
\right)_{j \in \bb{\mathbf{t}}}
\right)
\right)_{i\in\bb{\mathbf{r}}}
\right)
\right\rrbracket_{t}
\tag{7}
\\
&=
\left\llbracket
\tau^{\mathbf{Pth}_{\boldsymbol{\mathcal{B}}^{(2)}}^{(0,2)}} \left(
\left(
Q_{i}^{\mathbf{Pth}_{\boldsymbol{\mathcal{B}}^{(2)}}^{(0,2)}} \left(
\left(
\mathfrak{Q}_{j}^{(2)}
\right)_{j \in \bb{\mathbf{t}}}
\right)
\right)_{i\in\bb{\mathbf{r}}}
\right)
\right.
\\
&\hspace{3.1cm}
\left.
\circ_{t}^{0\mathbf{Pth}_{\boldsymbol{\mathcal{B}}^{(2)}}}
\tau^{\mathbf{Pth}_{\boldsymbol{\mathcal{B}}^{(2)}}^{(0,2)}} \left(
\left(
Q_{i}^{\mathbf{Pth}_{\boldsymbol{\mathcal{B}}^{(2)}}^{(0,2)}} \left(
\left(
\mathfrak{P}_{j}^{(2)}
\right)_{j \in \bb{\mathbf{t}}}
\right)
\right)_{i\in\bb{\mathbf{r}}}
\right)
\right\rrbracket_{t}
\tag{8}
\\
&=
\left\llbracket
\tau^{\mathbf{O}_{\mathbf{t}}\left(\mathbf{Pth}_{\boldsymbol{\mathcal{B}}^{(2)}}^{(0,2)}\right)} \left(
\left(
Q_{i}^{\mathbf{Pth}_{\boldsymbol{\mathcal{B}}^{(2)}}^{(0,2)}}
\right)_{i \in \bb{\mathbf{r}}}
\right)
\left(
\left(
\mathfrak{Q}_{j}^{(2)}
\right)_{j\in\bb{\mathbf{t}}}
\right)
\right.
\\
&\hspace{2.5cm}
\left.
\circ_{t}^{0\mathbf{Pth}_{\boldsymbol{\mathcal{B}}^{(2)}}}
\tau^{\mathbf{O}_{\mathbf{t}}\left(\mathbf{Pth}_{\boldsymbol{\mathcal{B}}^{(2)}}^{(0,2)}\right)} \left(
\left(
Q_{i}^{\mathbf{Pth}_{\boldsymbol{\mathcal{B}}^{(2)}}^{(0,2)}}
\right)_{i \in \bb{\mathbf{r}}}
\right)
\left(
\left(
\mathfrak{P}_{j}^{(2)}
\right)_{j\in\bb{\mathbf{t}}}
\right)
\right\rrbracket_{t}
\tag{9}
\\
&=
\left\llbracket
P^{\mathbf{Pth}_{\boldsymbol{\mathcal{B}}^{(2)}}^{(0,2)}}
\left(
\left(
\mathfrak{Q}_{j}^{(2)}
\right)_{j\in\bb{\mathbf{t}}}
\right)
\circ_{t}^{0\mathbf{Pth}_{\boldsymbol{\mathcal{B}}^{(2)}}}
P^{\mathbf{Pth}_{\boldsymbol{\mathcal{B}}^{(2)}}^{(0,2)}}
\left(
\left(
\mathfrak{P}_{j}^{(2)}
\right)_{j\in\bb{\mathbf{t}}}
\right)
\right\rrbracket_{t}.
\tag{10}
\end{align*}
\end{flushleft}

The first equality unravels the definition of $P$;
the second equality unravels the definition of the interpretation of the operation symbol $\tau$ in the product $\Lambda$-algebra $\mathbf{O}_{\mathbf{t}}\left(\mathbf{Pth}_{\boldsymbol{\mathcal{B}}^{(2)}}^{(0,2)}\right)$ introduced in Definition~\ref{wdop};
the third equality follows from the definition of the interpretation of the operation symbol $\tau$ in the quotient $\Lambda$-algebra $\left\llbracket\mathbf{Pth}_{\boldsymbol{\mathcal{B}}^{(2)}}^{(0,2)}\right\rrbracket$ introduced in Proposition~\ref{PDVDCatAlg};
the fourth equality follows by the inductive hypothesis;
the fifth equality follows from the definition of the interpretation of the operation symbol $\circ_{r_{i}}^{0}$ in the quotient $\Lambda^{\boldsymbol{\mathcal{B}}^{(2)}}$-algebra $\llbracket\mathbf{Pth}_{\boldsymbol{\mathcal{B}}^{(2)}}\rrbracket$ introduced in Proposition~\ref{PDVDCatAlg};
the sixth equality follows from Proposition~\ref{PDVVarA8};
the seventh equality recovers the definition of the interpretation of the operation symbol $\tau$ in the quotient $\Lambda$-algebra $\llbracket\mathbf{Pth}_{\boldsymbol{\mathcal{B}}^{(2)}}^{(0,2)}\rrbracket$ introduced in Proposition~\ref{PDVDCatAlg};
the eighth equality recovers the definition of the interpretation of the operation symbol $\circ_{t}^{0}$ in the quotient $\Lambda^{\boldsymbol{\mathcal{B}}^{(2)}}$-algebra $\llbracket\mathbf{Pth}_{\boldsymbol{\mathcal{B}}^{(2)}}\rrbracket$ introduced in Proposition~\ref{PDVDCatAlg};
the ninth equality recovers the definition of the interpretation of the operation symbol $\tau$ in the product $\Lambda$-algebra $\mathbf{O}_{\mathbf{t}}\left(\mathbf{Pth}_{\boldsymbol{\mathcal{B}}^{(2)}}^{(0,2)}\right)$ introduced in Definition~\ref{wdop};
finally, the last equality recovers the definition of $P$.
\end{proof}

\begin{proposition}
\label{PDQPthBVarB8}
Let $(\mathbf{s}, s)$ an element of $S^{\star} \times S$, $\sigma$ an operation symbol in $\Sigma_{\mathbf{s}, s}$ and $\left( \llbracket\mathfrak{P}_{j}^{(2)}\rrbracket_{\varphi(s_{j})} \right)_{j\in\bb{\mathbf{s}}}$ and $\left( \llbracket\mathfrak{Q}_{j}^{(2)}\rrbracket_{\varphi(s_{j})} \right)_{j\in\bb{\mathbf{s}}}$ be two family of path classes in $\llbracket\mathrm{Pth}_{\boldsymbol{\mathcal{B}}^{(2)}}\rrbracket_{\varphi^{\star}(\mathbf{s})}$ such that, for every $j \in \bb{\mathbf{s}}$, 
$$
\mathrm{sc}_{s_{j}}^{1\llbracket\mathbf{Pth}_{\boldsymbol{\mathcal{B}}^{(2)}}^{\mathbf{f}^{(2)}}\rrbracket} \left(
\left\llbracket
\mathfrak{Q}_{j}^{(2)}
\right\rrbracket_{s_{j}}
\right)
=
\mathrm{tg}_{s_{j}}^{1\llbracket\mathbf{Pth}_{\boldsymbol{\mathcal{B}}^{(1)}}^{\mathbf{f}^{(1)}}\rrbracket} \left(
\left\llbracket
\mathfrak{P}_{j}^{(2)}
\right\rrbracket_{s_{j}}
\right).
$$
Then the following equality holds
\begin{multline*}
\sigma^{\llbracket\mathbf{Pth}_{\boldsymbol{\mathcal{B}}^{(2)}}^{\mathbf{f}^{(2)}}\rrbracket} \left(
\left(
\left\llbracket
\mathfrak{Q}_{j}^{(2)}
\right\rrbracket_{s_{j}}
\circ_{s_{j}}^{1\llbracket\mathbf{Pth}_{\boldsymbol{\mathcal{B}}^{(2)}}^{\mathbf{f}^{(2)}}\rrbracket}
\left\llbracket
\mathfrak{P}_{j}^{(2)}
\right\rrbracket_{s_{j}}
\right)_{j \in \bb{\mathbf{s}}}
\right)
\\
=
\sigma^{\llbracket\mathbf{Pth}_{\boldsymbol{\mathcal{B}}^{(2)}}^{\mathbf{f}^{(2)}}\rrbracket} \left(
\left(
\left\llbracket
\mathfrak{Q}_{j}^{(2)}
\right\rrbracket_{s_{j}}
\right)_{j \in \bb{\mathbf{s}}}
\right)
\circ_{s}^{1\llbracket\mathbf{Pth}_{\boldsymbol{\mathcal{B}}^{(2)}}^{\mathbf{f}^{(2)}}\rrbracket}
\sigma^{\llbracket\mathbf{Pth}_{\boldsymbol{\mathcal{B}}^{(2)}}^{\mathbf{f}^{(2)}}\rrbracket} \left(
\left(
\left\llbracket
\mathfrak{P}_{j}^{(2)}
\right\rrbracket_{s_{j}}
\right)_{j \in \bb{\mathbf{s}}}
\right).
\end{multline*}
\end{proposition}

\begin{proof}
Let $(\mathbf{s}, s)$ an element of $S^{\star} \times S$, $\sigma$ an operation symbol in $\Sigma_{\mathbf{s}, s}$ and $( \llbracket\mathfrak{P}^{(2)}_{j}\rrbracket_{\varphi(s_{j})} )_{j\in\bb{\mathbf{s}}}$ and $( \llbracket\mathfrak{Q}^{(2)}_{j}\rrbracket_{\varphi(s_{j})} )_{j\in\bb{\mathbf{s}}}$ be two family of path classes in $\llbracket\mathrm{Pth}_{\boldsymbol{\mathcal{B}}^{(2)}}\rrbracket_{\varphi^{\star}(\mathbf{s})}$ such that, for every $j \in \bb{\mathbf{s}}$, 
$$
\mathrm{sc}_{s_{j}}^{1\llbracket\mathbf{Pth}_{\boldsymbol{\mathcal{B}}^{(2)}}^{\mathbf{f}^{(2)}}\rrbracket} \left(
\left\llbracket
\mathfrak{Q}_{j}^{(2)}
\right\rrbracket_{s_{j}}
\right)
=
\mathrm{tg}_{s_{j}}^{1\llbracket\mathbf{Pth}_{\boldsymbol{\mathcal{B}}^{(1)}}^{\mathbf{f}^{(1)}}\rrbracket} \left(
\left\llbracket
\mathfrak{P}_{j}^{(2)}
\right\rrbracket_{s_{j}}
\right).
$$
Unraveling the definition of the interpretation of the operation symbols $\sigma$ and $\circ_{s}^{1}$ in the $\Sigma^{\boldsymbol{\mathcal{A}}^{(2)}}$-algebras $\llbracket\mathbf{Pth}_{\boldsymbol{\mathcal{B}}^{(2)}}^{\mathbf{f}^{(2)}}\rrbracket$ and $\mathbf{Pth}_{\boldsymbol{\mathcal{B}}^{(2)}}^{\mathbf{f}^{(2)}}$ introduced in Propositions~\ref{PDQPthBDCatAlg} and \ref{PDPthBCatAlg}, respectively, the desired equality is equivalent to
\begin{multline*}
\left\llbracket
\sigma^{\mathbf{c}_{\mathfrak{d}}^{\ast}(\mathbf{Pth}_{\boldsymbol{\mathcal{B}}^{(2)}}^{(0,2)})} \left(
\left(
\mathfrak{Q}_{j}^{(2)}
\circ_{\varphi(s_{j})}^{1\mathbf{Pth}_{\boldsymbol{\mathcal{B}}^{(2)}}}
\mathfrak{P}_{j}^{(2)}
\right)_{j\in\bb{\mathbf{s}}}
\right)
\right\rrbracket_{\varphi(s)}
\\
=
\left\llbracket
\sigma^{\mathbf{c}_{\mathfrak{d}}^{\ast}(\mathbf{Pth}_{\boldsymbol{\mathcal{B}}^{(2)}}^{(0,2)})} \left(
\left(
\mathfrak{Q}_{j}^{(2)}
\right)_{j \in \bb{\mathbf{s}}}
\right)
\circ_{\varphi(s)}^{1\mathbf{Pth}_{\boldsymbol{\mathcal{B}}^{(2)}}}
\sigma^{\mathbf{c}_{\mathfrak{d}}^{\ast}(\mathbf{Pth}_{\boldsymbol{\mathcal{B}}^{(2)}}^{(0,2)})} \left(
\left(
\mathfrak{P}_{j}^{(2)}
\right)_{j \in \bb{\mathbf{s}}}
\right)
\right\rrbracket_{\varphi(s)}.
\tag{\dag}
\end{multline*}

Let us recall that, the interpretation of the operation symbol $\sigma$ in the $\Sigma$-algebra $\mathbf{c}_{\mathfrak{d}}^{\ast}(\mathbf{Pth}_{\boldsymbol{\mathcal{B}}^{(2)}}^{(0,2)})$ is given by the derived operation on $\mathbf{Pth}_{\boldsymbol{\mathcal{B}}^{(2)}}^{(0,2)}$ determined by $c_{\mathbf{s}, s}(\sigma)$ which following Remark~\ref{Rdop}, we denote by $c(\sigma)^{\mathbf{Pth}_{\boldsymbol{\mathcal{B}}^{(2)}}^{(0,2)}}$. For the details about derived operations, see Definition~\ref{wdop}.

Thus, we prove that, for every $(\mathbf{t}, t)$ in $T^{\star} \times T$, every $P$ in $T_{\Lambda}(\vs\mathbf{t})_{t}$ and every pair of families $(\mathfrak{P}_{j}^{(2)})_{j \in \bb{\mathbf{t}}}$ and $(\mathfrak{Q}_{j}^{(2)})_{j \in \bb{\mathbf{t}}}$ in $\mathrm{Pth}_{\boldsymbol{\mathcal{B}}^{(2)}, \mathbf{t}}$,
\begin{multline*}
\left\llbracket
P^{\mathbf{Pth}_{\boldsymbol{\mathcal{B}}^{(2)}}^{(0,2)}} \left(
\left(
\mathfrak{Q}_{j}^{(2)}
\circ_{t_{j}}^{1\mathbf{Pth}_{\boldsymbol{\mathcal{B}}^{(2)}}}
\mathfrak{P}_{j}^{(2)}
\right)_{j\in\bb{\mathbf{t}}}
\right)
\right\rrbracket_{t}
\\
=
\left\llbracket
P^{\mathbf{Pth}_{\boldsymbol{\mathcal{B}}^{(2)}}^{(0,2)}} \left(
\left(
\mathfrak{Q}_{j}^{(2)}
\right)_{j \in \bb{\mathbf{t}}}
\right)
\circ_{t}^{1\mathbf{Pth}_{\boldsymbol{\mathcal{B}}^{(2)}}}
P^{\mathbf{Pth}_{\boldsymbol{\mathcal{B}}^{(2)}}^{(0,2)}} \left(
\left(
\mathfrak{P}_{j}^{(2)}
\right)_{j \in \bb{\mathbf{t}}}
\right)
\right\rrbracket_{t}.
\end{multline*}
Therefore, Equation~(\dag) follows as a consequence.

By algebraic induction on the complexity of $P$.

{\sffamily Base step of the induction}

If $P$ is a variable $v_{i}^{t}$, thus, $i \in \bb{t}$ and $t_{i} = t$, then following Remark~\ref{Rdop}
$$
(v_{i}^{t})^{\mathbf{Pth}_{\boldsymbol{\mathcal{B}}^{(2)}}^{(0,2)}}
=
\mathrm{d}_{\mathbf{t}, t}^{\mathbf{Pth}_{\boldsymbol{\mathcal{B}}^{(2)}}^{(0,2)}} (v_{i}^{t})
=
\mathrm{d}_{\mathbf{t}, t}^{\mathbf{Pth}_{\boldsymbol{\mathcal{B}}^{(2)}}^{(0,2)}} ( \eta^{\vs\mathbf{t}} (v_{i}^{t})) 
=
\mathrm{p}_{\mathbf{t}, t}^{\mathbf{Pth}_{\boldsymbol{\mathcal{B}}^{(2)}}^{(0,2)}} (v_{i}^{t})
=
\mathrm{pr}_{\mathbf{t}, i}^{\mathrm{Pth}_{\boldsymbol{\mathcal{B}}^{(2)}}}.
$$

The following chain of equalities holds
\allowdisplaybreaks
\begin{flushleft}
$
P^{\mathbf{Pth}_{\boldsymbol{\mathcal{B}}^{(2)}}^{(0,2)}} \left(
\left(
\mathfrak{Q}_{j}^{(2)}
\circ_{t_{j}}^{1\mathbf{Pth}_{\boldsymbol{\mathcal{B}}^{(2)}}}
\mathfrak{P}_{j}^{(2)}
\right)_{j\in\bb{\mathbf{t}}}
\right)
$
\begin{align*}
&=
(v_{i}^{t})^{\mathbf{Pth}_{\boldsymbol{\mathcal{B}}^{(2)}}^{(0,2)}} \left(
\left(
\mathfrak{Q}_{j}^{(2)}
\circ_{t_{j}}^{1\mathbf{Pth}_{\boldsymbol{\mathcal{B}}^{(2)}}}
\mathfrak{P}_{j}^{(2)}
\right)_{j\in\bb{\mathbf{t}}}
\right)
\tag{1}
\\
&=
\mathrm{pr}_{\mathbf{t}, i}^{\mathrm{Pth}_{\boldsymbol{\mathcal{B}}^{(2)}}} \left(
\left(
\mathfrak{Q}_{j}^{(2)}
\circ_{t_{j}}^{1\mathbf{Pth}_{\boldsymbol{\mathcal{B}}^{(2)}}}
\mathfrak{P}_{j}^{(2)}
\right)_{j\in\bb{\mathbf{t}}}
\right)
\tag{2}
\\
&=
\mathfrak{Q}_{i}^{(2)}
\circ_{t_{i}}^{1\mathbf{Pth}_{\boldsymbol{\mathcal{B}}^{(2)}}}
\mathfrak{P}_{i}^{(2)}
\tag{3}
\\
&=
\mathfrak{Q}_{i}^{(2)}
\circ_{t}^{1\mathbf{Pth}_{\boldsymbol{\mathcal{B}}^{(2)}}}
\mathfrak{P}_{i}^{(2)}
\tag{4}
\\
&=
\mathrm{pr}_{\mathbf{t}, i}^{\mathrm{Pth}_{\boldsymbol{\mathcal{B}}^{(2)}}} \left(
\left(
\mathfrak{Q}_{j}^{(2)}
\right)_{j \in \bb{\mathbf{t}}}
\right)
\circ_{t}^{1\mathbf{Pth}_{\boldsymbol{\mathcal{B}}^{(2)}}}
\mathrm{pr}_{\mathbf{t}, i}^{\mathrm{Pth}_{\boldsymbol{\mathcal{B}}^{(2)}}} \left(
\left(
\mathfrak{P}_{j}^{(2)}
\right)_{j\in\bb{\mathbf{t}}}
\right)
\tag{5}
\\
&=
(v_{i}^{t})^{\mathbf{Pth}_{\boldsymbol{\mathcal{B}}^{(2)}}^{(0,2)}} \left(
\left(
\mathfrak{Q}_{j}^{(2)}
\right)_{j \in \bb{\mathbf{t}}}
\right)
\circ_{t}^{1\mathbf{Pth}_{\boldsymbol{\mathcal{B}}^{(2)}}}
(v_{i}^{t})^{\mathbf{Pth}_{\boldsymbol{\mathcal{B}}^{(2)}}^{(0,2)}} \left(
\left(
\mathfrak{P}_{j}^{(2)}
\right)_{j\in\bb{\mathbf{t}}}
\right)
\tag{6}
\\
&=
P^{\mathbf{Pth}_{\boldsymbol{\mathcal{B}}^{(2)}}^{(0,2)}} \left(
\left(
\mathfrak{Q}_{j}^{(2)}
\right)_{j \in \bb{\mathbf{t}}}
\right)
\circ_{t}^{1\mathbf{Pth}_{\boldsymbol{\mathcal{B}}^{(2)}}}
P^{\mathbf{Pth}_{\boldsymbol{\mathcal{B}}^{(2)}}^{(0,2)}} \left(
\left(
\mathfrak{P}_{j}^{(2)}
\right)_{j\in\bb{\mathbf{t}}}
\right).
\tag{7}
\end{align*}
\end{flushleft}

The first equality unravels the definition of $P$;
the second equality follows from the fact that $ (v_{i}^{t})^{\mathbf{Pth}_{\boldsymbol{\mathcal{B}}^{(2)}}^{(0,2)}} = \mathrm{pr}_{i}^{\mathrm{Pth}_{\boldsymbol{\mathcal{B}}^{(2)}, \mathbf{t}}} $;
the third equality unravels the definition of the $\mathbf{t}$-ary, $i$-th canonical projection from $\mathrm{Pth}_{\boldsymbol{\mathcal{B}}^{(2)}, \mathbf{t}}$ to $\mathrm{Pth}_{\boldsymbol{\mathcal{B}}^{(2)}, t_{i}}$;
the fourth equality follows from the fact that $t_{i} = t$;
the fifth equality recovers the definition of the $\mathbf{t}$-ary, $i$-th canonical projection from $\mathrm{Pth}_{\boldsymbol{\mathcal{B}}^{(2)}, \mathbf{t}}$ to $\mathrm{Pth}_{\boldsymbol{\mathcal{B}}^{(2)}, t_{i}}$;
the sixth equality follows from the fact that $ (v_{i}^{t})^{\mathbf{Pth}_{\boldsymbol{\mathcal{B}}^{(2)}}^{(0,2)}} = \mathrm{pr}_{i}^{\mathrm{Pth}_{\boldsymbol{\mathcal{B}}^{(2)}, \mathbf{t}}} $;
finally, the last equality recovers the definition of $P$.

Therefore, its $\llbracket\cdot\rrbracket$-classes coincide, i.e., equality (\dag) follows.

If $P$ is a constant symbol $\tau$ in $\Lambda_{\lambda, t}$, then the following chain of equalities holds
\allowdisplaybreaks
\begin{flushleft}
$
P^{\mathbf{Pth}_{\boldsymbol{\mathcal{B}}^{(2)}}^{(0,2)}} \left(
\left(
\mathfrak{Q}_{j}^{(2)}
\circ_{t_{j}}^{1\mathbf{Pth}_{\boldsymbol{\mathcal{B}}^{(2)}}}
\mathfrak{P}_{j}^{(2)}
\right)_{j\in\bb{\mathbf{t}}}
\right)
$
\begin{align*}
&=
\tau^{\mathbf{Pth}_{\boldsymbol{\mathcal{B}}^{(2)}}^{(0,2)}}
\tag{1}
\\
&=
\tau^{\mathbf{Pth}_{\boldsymbol{\mathcal{B}}^{(2)}}^{(0,2)}}
\circ_{t}^{1\mathbf{Pth}_{\boldsymbol{\mathcal{B}}^{(2)}}}
\tau^{\mathbf{Pth}_{\boldsymbol{\mathcal{B}}^{(2)}}^{(0,2)}}
\tag{2}
\\
&=
P^{\mathbf{Pth}_{\boldsymbol{\mathcal{B}}^{(2)}}^{(0,2)}} \left(
\left(
\mathfrak{Q}_{j}^{(2)}
\right)_{j \in \bb{\mathbf{t}}}
\right)
\circ_{t}^{1\mathbf{Pth}_{\boldsymbol{\mathcal{B}}^{(2)}}}
P^{\mathbf{Pth}_{\boldsymbol{\mathcal{B}}^{(2)}}^{(0,2)}} \left(
\left(
\mathfrak{P}_{j}^{(2)}
\right)_{j\in\bb{\mathbf{t}}}
\right).
\tag{3}
\end{align*}
\end{flushleft}

The first equality unravels the definition of $P$;
the second equality follows from the fact that, according to Proposition~\ref{PDPthCatAlg}, the interpretation of the $\tau$ operation in the many-sorted partial $\Lambda^{\boldsymbol{\mathcal{B}}^{(2)}}$-algebra $\mathbf{Pth}_{\boldsymbol{\mathcal{B}}^{(2)}}$ is given by the $(2,[1])$-identity path on $[\tau^{\mathbf{PT}_{\boldsymbol{\mathcal{B}}^{(1)}}}]$, i.e., $\tau^{\mathbf{Pth}_{\boldsymbol{\mathcal{B}}^{(2)}}} = \mathrm{ip}_{t}^{(2,[1])\sharp} \left([\tau^{\mathbf{PT}_{\boldsymbol{\mathcal{B}}^{(1)}}}]\right)$ and from the fact that, according to Proposition~\ref{PDPthComp}, the $(2,[1])$-identity paths are idempotent for the $1$-composition;
finally, the last equality recovers the definition of $P$.

Therefore, its $\llbracket\cdot\rrbracket$-classes coincide, i.e., equality (\dag) follows.

{\sffamily Inductive step of the induction}

If $P$ is a term $\tau ((Q_{i})_{i \in \bb{\mathbf{r}}})$ for $\tau$ an operation symbol in $\lambda_{\mathbf{r}, t}$ and $(Q_{i})_{i \in \bb{\mathbf{r}}}$ a family of terms in $\T_{\Lambda}(\vs\mathbf{t})_{\mathbf{r}}$, then following Remark~\ref{Rdop}
\begin{align*}
P^{\mathbf{Pth}_{\boldsymbol{\mathcal{B}}^{(2)}}^{(0,2)}}
&=
\mathrm{d}_{\mathbf{t}, t}^{\mathbf{Pth}_{\boldsymbol{\mathcal{B}}^{(2)}}^{(0,2)}} (P)
\\
&=
\mathrm{d}_{\mathbf{t}, t}^{\mathbf{Pth}_{\boldsymbol{\mathcal{B}}^{(2)}}^{(0,2)}} \left(\tau \left(\left(Q_{i}\right)_{i \in \bb{\mathbf{r}}}\right)\right) 
\\
&=
\tau^{\mathbf{O}_{\mathbf{t}}(\mathbf{Pth}_{\boldsymbol{\mathcal{B}}^{(2)}}^{(0,2)})} \left(
\left(
\mathrm{d}_{\mathbf{t}, r_{i}}^{\mathbf{Pth}_{\boldsymbol{\mathcal{B}}^{(2)}}^{(0,2)}}
\left(
Q_{i}
\right)
\right)_{i \in \bb{\mathbf{r}}}
\right)
\\
&=
\tau^{\mathbf{O}_{\mathbf{t}}(\mathbf{Pth}_{\boldsymbol{\mathcal{B}}^{(2)}}^{(0,2)})} \left(
\left(
Q_{i}^{\mathbf{Pth}_{\boldsymbol{\mathcal{B}}^{(2)}}^{(0,2)}}
\right)_{i \in \bb{\mathbf{r}}}
\right).
\end{align*}

The following chain of equalities holds
\allowdisplaybreaks
\begin{flushleft}
$
\left\llbracket
P^{\mathbf{Pth}_{\boldsymbol{\mathcal{B}}^{(2)}}^{(0,2)}} \left(
\left(
\mathfrak{Q}_{j}^{(2)}
\circ_{t_{j}}^{1\mathbf{Pth}_{\boldsymbol{\mathcal{B}}^{(2)}}}
\mathfrak{P}_{j}^{(2)}
\right)_{j\in\bb{\mathbf{t}}}
\right)
\right\rrbracket_{t}
$
\begin{align*}
&=
\left\llbracket
\tau^{\mathbf{O}_{\mathbf{t}}(\mathbf{Pth}_{\boldsymbol{\mathcal{B}}^{(2)}}^{(0,2)})} \left(
\left(
Q_{i}^{\mathbf{Pth}_{\boldsymbol{\mathcal{B}}^{(2)}}^{(0,2)}}
\right)_{i \in \bb{\mathbf{r}}}
\right)
\left(
\left(
\mathfrak{Q}_{j}^{(2)}
\circ_{t_{j}}^{1\mathbf{Pth}_{\boldsymbol{\mathcal{B}}^{(2)}}}
\mathfrak{P}_{j}^{(2)}
\right)_{j\in\bb{\mathbf{t}}}
\right)
\right\rrbracket_{t}
\tag{1}
\\
&=
\left\llbracket
\tau^{\mathbf{Pth}_{\boldsymbol{\mathcal{B}}^{(2)}}^{(0,2)}} \left(
\left(
Q_{i}^{\mathbf{Pth}_{\boldsymbol{\mathcal{B}}^{(2)}}^{(0,2)}}
\left(
\left(
\mathfrak{Q}_{j}^{(2)}
\circ_{t_{j}}^{1\mathbf{Pth}_{\boldsymbol{\mathcal{B}}^{(2)}}}
\mathfrak{P}_{j}^{(2)}
\right)_{j\in\bb{\mathbf{t}}}
\right)
\right)_{i \in \bb{\mathbf{r}}}
\right)
\right\rrbracket_{t}
\tag{2}
\\
&=
\tau^{\llbracket\mathbf{Pth}_{\boldsymbol{\mathcal{B}}^{(2)}}^{(0,2)}\rrbracket} \left(
\left(
\left\llbracket
Q_{i}^{\mathbf{Pth}_{\boldsymbol{\mathcal{B}}^{(2)}}^{(0,2)}}
\left(
\left(
\mathfrak{Q}_{j}^{(2)}
\circ_{t_{j}}^{1\mathbf{Pth}_{\boldsymbol{\mathcal{B}}^{(2)}}}
\mathfrak{P}_{j}^{(2)}
\right)_{j\in\bb{\mathbf{t}}}
\right)
\right\rrbracket_{r_{i}}
\right)_{i\in\bb{\mathbf{r}}}
\right)
\tag{3}
\\
&=
\resizebox{0.89\textwidth}{!}{%
$
\tau^{\llbracket\mathbf{Pth}_{\boldsymbol{\mathcal{B}}^{(2)}}^{(0,2)}\rrbracket} \left(
\left(
\left\llbracket
Q_{i}^{\mathbf{Pth}_{\boldsymbol{\mathcal{B}}^{(2)}}^{(0,2)}} \left(
\left(
\mathfrak{Q}_{j}^{(2)}
\right)_{j \in \bb{\mathbf{t}}}
\right)
\circ_{r_{i}}^{1\mathbf{Pth}_{\boldsymbol{\mathcal{B}}^{(2)}}}
Q_{i}^{\mathbf{Pth}_{\boldsymbol{\mathcal{B}}^{(2)}}^{(0,2)}} \left(
\left(
\mathfrak{P}_{j}^{(2)}
\right)_{j \in \bb{\mathbf{t}}}
\right)
\right\rrbracket_{r_{i}}
\right)_{i\in\bb{\mathbf{r}}}
\right)
$}
\tag{4}
\\
&=
\resizebox{0.89\textwidth}{!}{%
$
\tau^{\llbracket\mathbf{Pth}_{\boldsymbol{\mathcal{B}}^{(2)}}^{(0,2)}\rrbracket} \left(
\left(
\left\llbracket
Q_{i}^{\mathbf{Pth}_{\boldsymbol{\mathcal{B}}^{(2)}}^{(0,2)}} \left(
\left(
\mathfrak{Q}_{j}^{(2)}
\right)_{j \in \bb{\mathbf{t}}}
\right)
\right\rrbracket_{r_{i}}
\circ_{r_{i}}^{1\left\llbracket\mathbf{Pth}_{\boldsymbol{\mathcal{B}}^{(2)}}\right\rrbracket}
\left\llbracket
Q_{i}^{\mathbf{Pth}_{\boldsymbol{\mathcal{B}}^{(2)}}^{(0,2)}} \left(
\left(
\mathfrak{P}_{j}^{(2)}
\right)_{j \in \bb{\mathbf{t}}}
\right)
\right\rrbracket_{r_{i}}
\right)_{i\in\bb{\mathbf{r}}}
\right)
$}
\tag{5}
\\
&=
\tau^{\llbracket\mathbf{Pth}_{\boldsymbol{\mathcal{B}}^{(2)}}^{(0,2)}\rrbracket} \left(
\left(
\left\llbracket
Q_{i}^{\mathbf{Pth}_{\boldsymbol{\mathcal{B}}^{(2)}}^{(0,2)}} \left(
\left(
\mathfrak{Q}_{j}^{(2)}
\right)_{j \in \bb{\mathbf{t}}}
\right)
\right\rrbracket_{r_{i}}
\right)_{i\in\bb{\mathbf{r}}}
\right)
\\
&\hspace{2cm}
\circ_{t}^{1\llbracket\mathbf{Pth}_{\boldsymbol{\mathcal{B}}^{(2)}}\rrbracket}
\tau^{\llbracket\mathbf{Pth}_{\boldsymbol{\mathcal{B}}^{(2)}}^{(0,2)}\rrbracket} \left(
\left(
\left\llbracket
Q_{i}^{\mathbf{Pth}_{\boldsymbol{\mathcal{B}}^{(2)}}^{(0,2)}} \left(
\left(
\mathfrak{P}_{j}^{(2)}
\right)_{j \in \bb{\mathbf{t}}}
\right)
\right\rrbracket_{r_{i}}
\right)_{i\in\bb{\mathbf{r}}}
\right)
\tag{6}
\\
&=
\left\llbracket
\tau^{\mathbf{Pth}_{\boldsymbol{\mathcal{B}}^{(2)}}^{(0,2)}} \left(
\left(
Q_{i}^{\mathbf{Pth}_{\boldsymbol{\mathcal{B}}^{(2)}}^{(0,2)}} \left(
\left(
\mathfrak{Q}_{j}^{(2)}
\right)_{j \in \bb{\mathbf{t}}}
\right)
\right)_{i\in\bb{\mathbf{r}}}
\right)
\right\rrbracket_{t}
\\
&\hspace{2.5cm}
\circ_{t}^{1\llbracket\mathbf{Pth}_{\boldsymbol{\mathcal{B}}^{(2)}}\rrbracket}
\left\llbracket
\tau^{\mathbf{Pth}_{\boldsymbol{\mathcal{B}}^{(2)}}^{(0,2)}} \left(
\left(
Q_{i}^{\mathbf{Pth}_{\boldsymbol{\mathcal{B}}^{(2)}}^{(0,2)}} \left(
\left(
\mathfrak{P}_{j}^{(2)}
\right)_{j \in \bb{\mathbf{t}}}
\right)
\right)_{i\in\bb{\mathbf{r}}}
\right)
\right\rrbracket_{t}
\tag{7}
\\
&=
\left\llbracket
\tau^{\mathbf{Pth}_{\boldsymbol{\mathcal{B}}^{(2)}}^{(0,2)}} \left(
\left(
Q_{i}^{\mathbf{Pth}_{\boldsymbol{\mathcal{B}}^{(2)}}^{(0,2)}} \left(
\left(
\mathfrak{Q}_{j}^{(2)}
\right)_{j \in \bb{\mathbf{t}}}
\right)
\right)_{i\in\bb{\mathbf{r}}}
\right)
\right.
\\
&\hspace{3.1cm}
\left.
\circ_{t}^{1\mathbf{Pth}_{\boldsymbol{\mathcal{B}}^{(2)}}}
\tau^{\mathbf{Pth}_{\boldsymbol{\mathcal{B}}^{(2)}}^{(0,2)}} \left(
\left(
Q_{i}^{\mathbf{Pth}_{\boldsymbol{\mathcal{B}}^{(2)}}^{(0,2)}} \left(
\left(
\mathfrak{P}_{j}^{(2)}
\right)_{j \in \bb{\mathbf{t}}}
\right)
\right)_{i\in\bb{\mathbf{r}}}
\right)
\right\rrbracket_{t}
\tag{8}
\\
&=
\left\llbracket
\tau^{\mathbf{O}_{\mathbf{t}}(\mathbf{Pth}_{\boldsymbol{\mathcal{B}}^{(2)}}^{(0,2)})} \left(
\left(
Q_{i}^{\mathbf{Pth}_{\boldsymbol{\mathcal{B}}^{(2)}}^{(0,2)}}
\right)_{i \in \bb{\mathbf{r}}}
\right)
\left(
\left(
\mathfrak{Q}_{j}^{(2)}
\right)_{j\in\bb{\mathbf{t}}}
\right)
\right.
\\
&\hspace{2.5cm}
\left.
\circ_{t}^{1\mathbf{Pth}_{\boldsymbol{\mathcal{B}}^{(2)}}}
\tau^{\mathbf{O}_{\mathbf{t}}(\mathbf{Pth}_{\boldsymbol{\mathcal{B}}^{(2)}}^{(0,2)})} \left(
\left(
Q_{i}^{\mathbf{Pth}_{\boldsymbol{\mathcal{B}}^{(2)}}^{(0,2)}}
\right)_{i \in \bb{\mathbf{r}}}
\right)
\left(
\left(
\mathfrak{P}_{j}^{(2)}
\right)_{j\in\bb{\mathbf{t}}}
\right)
\right\rrbracket_{t}
\tag{9}
\\
&=
\left\llbracket
P^{\mathbf{Pth}_{\boldsymbol{\mathcal{B}}^{(2)}}^{(0,2)}}
\left(
\left(
\mathfrak{Q}_{j}^{(2)}
\right)_{j\in\bb{\mathbf{t}}}
\right)
\circ_{t}^{1\mathbf{Pth}_{\boldsymbol{\mathcal{B}}^{(2)}}}
P^{\mathbf{Pth}_{\boldsymbol{\mathcal{B}}^{(2)}}^{(0,2)}}
\left(
\left(
\mathfrak{P}_{j}^{(2)}
\right)_{j\in\bb{\mathbf{t}}}
\right)
\right\rrbracket_{t}.
\tag{10}
\end{align*}
\end{flushleft}

The first equality unravels the definition of $P$;
the second equality unravels the definition of the interpretation of the operation symbol $\tau$ in the product $\Lambda$-algebra $\mathbf{O}_{\mathbf{t}}(\mathbf{Pth}_{\boldsymbol{\mathcal{B}}^{(2)}}^{(0,2)})$ introduced in Definition~\ref{wdop};
the third equality follows from the definition of the interpretation of the operation symbol $\tau$ in the quotient $\Lambda$-algebra $\llbracket\mathbf{Pth}_{\boldsymbol{\mathcal{B}}^{(2)}}^{(0,2)}\rrbracket$ introduced in Proposition~\ref{PDVDCatAlg};
the fourth equality follows by the inductive hypothesis;
the fifth equality follows from the definition of the interpretation of the operation symbol $\circ_{r_{i}}^{1}$ in the quotient $\Lambda^{\boldsymbol{\mathcal{B}}^{(2)}}$-algebra $\llbracket\mathbf{Pth}_{\boldsymbol{\mathcal{B}}^{(2)}}\rrbracket$ introduced in Proposition~\ref{PDVDCatAlg};
the sixth equality follows from Proposition~\ref{PDVVarA8};
the seventh equality recovers the definition of the interpretation of the operation symbol $\tau$ in the quotient $\Lambda$-algebra $\llbracket\mathbf{Pth}_{\boldsymbol{\mathcal{B}}^{(2)}}^{(0,2)}\rrbracket$ introduced in Proposition~\ref{PDVDCatAlg};
the eighth equality recovers the definition of the interpretation of the operation symbol $\circ_{t}^{1}$ in the quotient $\Lambda^{\boldsymbol{\mathcal{B}}^{(2)}}$-algebra $\llbracket\mathbf{Pth}_{\boldsymbol{\mathcal{B}}^{(2)}}\rrbracket$ introduced in Proposition~\ref{PDVDCatAlg};
the ninth equality recovers the definition of the interpretation of the operation symbol $\tau$ in the product $\Lambda$-algebra $\mathbf{O}_{\mathbf{t}}(\mathbf{Pth}_{\boldsymbol{\mathcal{B}}^{(2)}}^{(0,2)})$ introduced in Definition~\ref{wdop};
finally, the last equality recovers the definition of $P$.
\end{proof}

We now work towards the existence of the second-order quotient path-extension mapping.

\begin{proposition}\label{PDHomPthExtKer}
Let $\mathbf{f}^{(2)} = (\varphi,c,(f^{(i)})_{i\in 3})$ be a second-order morphism from $\boldsymbol{\mathcal{A}}^{(2)}$ to $\boldsymbol{\mathcal{B}}^{(2)}$. Then
$$
\mathrm{pr}_{\boldsymbol{\mathcal{B}}^{(2)}, \varphi}^{\llbracket\cdot\rrbracket} \circ f^{(2)\flat}
\colon
\mathrm{Pth}_{\boldsymbol{\mathcal{A}}^{(2)}}
\mor
\left\llbracket\mathrm{Pth}_{\boldsymbol{\mathcal{B}}^{(2)}}\right\rrbracket_{\varphi}
$$
is a $\Sigma^{\boldsymbol{\mathcal{A}}^{(2)}}$-homomorphism from $\mathbf{Pth}_{\boldsymbol{\mathcal{A}}^{(2)}}$ to $\llbracket\mathbf{Pth}_{\boldsymbol{\mathcal{B}}^{(2)}}^{\mathbf{f}^{(2)}}\rrbracket_{\varphi}$ satisfying that
$$
\mathrm{Ker}(\mathrm{CH}^{(2)}_{\boldsymbol{\mathcal{A}}^{(2)}})
\vee
\Upsilon^{[1]}
\subseteq
\mathrm{Ker}(\mathrm{pr}_{\boldsymbol{\mathcal{B}}^{(2)}, \varphi}^{\llbracket\cdot\rrbracket} \circ f^{(2)\flat}).
$$
For further reference, we will call this composition mapping, which will become useful in further constructions, the \emph{second-order homomorphic path extension mapping} of $\mathbf{f}^{(2)}$. Occasionally, we will denote it by $f^{\llbracket2\rrbracket\flat} = \mathrm{pr}_{\boldsymbol{\mathcal{B}}^{(2)}, \varphi}^{\llbracket\cdot\rrbracket} \circ f^{(2)\flat}$.
\end{proposition}

\begin{proof}
We prove that $\mathrm{pr}_{\boldsymbol{\mathcal{B}}^{(2)}, \varphi}^{\llbracket\cdot\rrbracket} \circ f^{(2)\flat}$ is compatible with every operation symbol in $\Sigma^{\boldsymbol{\mathcal{A}}^{(2)}}$.

{\sffamily The mapping $\mathrm{pr}_{\boldsymbol{\mathcal{B}}^{(2)}, \varphi}^{\llbracket\cdot\rrbracket} \circ f^{(2)\flat}$ is a $\Sigma$-homomorphism.}

By Proposition~\ref{PDPthExtHom}, $f^{(2)\flat}$ is a $\Sigma$-homomorphism from $\mathbf{Pth}_{\boldsymbol{\mathcal{A}}^{2}}^{(0,2)}$ to $\mathbf{Pth}_{\boldsymbol{\mathcal{B}}^{(2)}}^{\mathbf{f}^{(2)}(0,2)}$. Also, by Proposition~\ref{PDPrBCatHom}, $\mathrm{pr}_{\boldsymbol{\mathcal{B}}^{(2)}, \varphi}^{\llbracket\cdot\rrbracket}$ is a $\Sigma$-homomorphism from $\mathbf{Pth}_{\boldsymbol{\mathcal{B}}^{(2)}}^{\mathbf{f}^{(2)}(0,2)}$ to $\llbracket\mathbf{Pth}_{\boldsymbol{\mathcal{B}}^{(2)}}^{\mathbf{f}^{(2)}(0,2)}\rrbracket$. Thus, the composition $\mathrm{pr}_{\boldsymbol{\mathcal{B}}^{(2)}, \varphi}^{\llbracket\cdot\rrbracket} \circ f^{(2)\flat}$ is a $\Sigma$-homomorphism from $\mathbf{Pth}_{\boldsymbol{\mathcal{A}}^{(1)}}^{(0,2)}$ to $\llbracket\mathbf{Pth}_{\boldsymbol{\mathcal{B}}^{(2)}}^{\mathbf{f}^{(2)}(0,2)}\rrbracket$.

Therefore, for every word $\mathbf{s}$ in $S^{\star}$, every sort $s$ in $S$, every operation symbol $\sigma \in \Sigma_{\mathbf{s}, s}$ and every family $(\mathfrak{P}_{j}^{(2)})_{j \in \bb{\mathbf{s}}}$ in $\mathrm{Pth}_{\boldsymbol{\mathcal{A}}^{(2)}, \mathbf{s}}$,
$$
\left\llbracket
f_{s}^{(2)\flat}\left(
\sigma^{\mathbf{Pth}_{\boldsymbol{\mathcal{A}}^{(2)}}}\left(\left(
\mathfrak{P}_{j}^{(2)}
\right)_{j\in\bb{\mathbf{s}}}\right)
\right)
\right\rrbracket_{\varphi(s)}
=
\sigma^{\llbracket\mathbf{Pth}_{\boldsymbol{\mathcal{B}}^{(2)}}^{\mathbf{f}^{(2)}}\rrbracket} \left(\left(
\left\llbracket
f_{s_{j}}^{(2)\flat}\left(
\mathfrak{P}_{j}^{(2)}
\right)
\right\rrbracket_{\varphi(s_{j})}
\right)_{j\in\bb{\mathbf{s}}}\right).
$$

{\sffamily The mapping $\mathrm{pr}_{\boldsymbol{\mathcal{B}}^{(2)}, \varphi}^{\llbracket\cdot\rrbracket} \circ f^{(2)\flat}$ is compatible with the first-order rewrite rules.}

Let $s$ be a sort in $S$ and let $\mathfrak{p}$ be a rewrite rule in $\mathcal{A}_{s}^{(1)}$. We have to show that
$$
\left\llbracket
f_{s}^{(2)\flat}\left(
\mathfrak{p}^{\mathbf{Pth}_{\boldsymbol{\mathcal{A}}^{(2)}}}
\right)
\right\rrbracket_{\varphi(s)}
=
\mathfrak{p}^{\llbracket\mathbf{Pth}_{\boldsymbol{\mathcal{B}}^{(2)}}^{\mathbf{f}^{(2)}}\rrbracket}
$$

The following chain of equalities holds
\allowdisplaybreaks
\begin{align*}
\left\llbracket
f_{s}^{(2)\flat}\left(
\mathfrak{p}^{\mathbf{Pth}_{\boldsymbol{\mathcal{A}}^{(2)}}}
\right)
\right\rrbracket_{\varphi(s)}
&=
\left\llbracket
f_{s}^{(2)\flat}\left(
\mathrm{ech}_{\boldsymbol{\mathcal{A}}^{(2)}, s}^{(2,\mathcal{A}^{(1)})}\left(
\mathfrak{p}
\right)
\right)
\right\rrbracket_{\varphi(s)}
\tag{1}
\\
&=
\mathfrak{p}^{\llbracket\mathbf{Pth}_{\boldsymbol{\mathcal{B}}^{(2)}}^{\mathbf{f}^{(2)}}\rrbracket}.
\tag{2}
\end{align*}

The first equality unravels the interpretation of the constant symbol $\mathfrak{p}$ in the partial $\Sigma^{\boldsymbol{\mathcal{A}}^{(2)}}$-algebra $\mathbf{Pth}_{\boldsymbol{\mathcal{A}}^{(2)}}$, introduced in Proposition~\ref{PDPthCatAlg};
the second equality recovers the interpretation of the constant symbol $\mathfrak{p}$ in the partial $\Sigma^{\boldsymbol{\mathcal{A}}^{(2)}}$-algebra $\llbracket\mathbf{Pth}_{\boldsymbol{\mathcal{B}}^{(2)}}^{\mathbf{f}^{(2)}}\rrbracket$, introduced in Proposition~\ref{PDQPthBDCatAlg}.

Hence, $\mathrm{pr}_{\boldsymbol{\mathcal{B}}^{(2)}, \varphi}^{\llbracket\cdot\rrbracket} \circ f^{(2)\flat}$ is compatible with the first-order rewrite rules.

{\sffamily The mapping $\mathrm{pr}_{\boldsymbol{\mathcal{B}}^{(2)}, \varphi}^{\llbracket\cdot\rrbracket} \circ f^{(2)\flat}$ is compatible with the $0$-source.}

Let $s$ be a sort in $S$ and let us consider the $0$-source operation symbol $\mathrm{sc}_{s}^{0}$ in $\Sigma^{\boldsymbol{\mathcal{A}}^{(2)}}_{s,s}$. Let $\mathfrak{P}^{(2)}$ be a second-order path in $\mathrm{Pth}_{\boldsymbol{\mathcal{A}}^{(2)}, s}$.

The following chain of equalities holds
\begin{flushleft}
$
\left\llbracket
f_{s}^{(2)\flat} \left(
\mathrm{sc}_{s}^{0\mathbf{Pth}_{\boldsymbol{\mathcal{A}}^{(2)}}} \left(
\mathfrak{P}^{(2)}
\right)
\right)
\right\rrbracket_{\varphi(s)}
$
\allowdisplaybreaks
\begin{align*}
&=
\left\llbracket
f_{s}^{(2)\flat} \left(
\mathrm{ip}_{\boldsymbol{\mathcal{A}}^{(2)}, s}^{(2,0)\sharp} \left(
\mathrm{sc}_{\boldsymbol{\mathcal{A}}^{(2)}, s}^{(0,2)} \left(
\mathfrak{P}^{(2)}
\right)
\right)
\right)
\right\rrbracket_{\varphi(s)}
\tag{1}
\\
&=
\left\llbracket
\mathrm{ip}_{\boldsymbol{\mathcal{B}}^{(2)}, \varphi(s)}^{(2,0)\sharp} \left(
f_{s}^{(0)\sharp} \left(
\mathrm{sc}_{\boldsymbol{\mathcal{A}}^{(2)}, s}^{(0,2)} \left(
\mathfrak{P}^{(2)}
\right)
\right)
\right)
\right\rrbracket_{\varphi(s)}
\tag{2}
\\
&=
\left\llbracket
\mathrm{ip}_{\boldsymbol{\mathcal{B}}^{(2)}, \varphi(s)}^{(2,0)\sharp} \left(
\mathrm{sc}_{\boldsymbol{\mathcal{B}}^{(2)}, \varphi(s)}^{(0,2)} \left(
f_{s}^{(2)\flat} \left(
\mathfrak{P}^{(2)}
\right)
\right)
\right)
\right\rrbracket_{\varphi(s)}
\tag{3}
\\
&=
\left\llbracket
\mathrm{sc}_{\varphi(s)}^{0\mathbf{Pth}_{\boldsymbol{\mathcal{B}}^{(2)}}} \left(
f_{s}^{(2)\flat} \left(
\mathfrak{P}^{(2)}
\right)
\right)
\right\rrbracket_{\varphi(s)}
\tag{4}
\\
&=
\left\llbracket
\mathrm{sc}_{s}^{0\mathbf{Pth}_{\boldsymbol{\mathcal{B}}^{(2)}}^{\mathbf{f}^{(2)}}} \left(
f_{s}^{(2)\flat} \left(
\mathfrak{P}^{(2)}
\right)
\right)
\right\rrbracket_{\varphi(s)}
\tag{5}
\\
&=
\mathrm{sc}_{s}^{0\llbracket\mathbf{Pth}_{\boldsymbol{\mathcal{B}}^{(2)}}^{\mathbf{f}^{(2)}}\rrbracket} \left(
\left\llbracket
f_{s}^{(2)\flat} \left(
\mathfrak{P}^{(2)}
\right)
\right\rrbracket_{\varphi(s)}
\right).
\tag{6}
\end{align*}
\end{flushleft}

The first equality unravels the interpretation of the operation symbol $\mathrm{sc}_{s}^{0}$ in the partial $\Sigma^{\boldsymbol{\mathcal{A}}^{(2)}}$-algebra $\mathbf{Pth}_{\boldsymbol{\mathcal{A}}^{(2)}}$, introduced in Proposition~\ref{PDPthDCatAlg}; 
the second equality follows from Proposition~\ref{PDPthExtDZIp}; 
the third equality follows from Proposition~\ref{PDPthExtDZScTg}; 
the fourth equality recovers the interpretation of the operation symbol $\mathrm{sc}_{\varphi(s)}^{0}$ in the partial $\Lambda^{\boldsymbol{\mathcal{B}}^{(2)}}$-algebra $\mathbf{Pth}_{\boldsymbol{\mathcal{B}}^{(2)}}$, according to Proposition~\ref{PDPthDCatAlg}; 
the fifth equality recovers the interpretation of the operation symbol $\mathrm{sc}_{s}^{0}$ in the partial $\Sigma^{\boldsymbol{\mathcal{A}}^{(2)}}$-algebra $\mathbf{Pth}_{\boldsymbol{\mathcal{B}}^{(2)}}^{\mathbf{f}^{(2)}}$, introduced in Proposition~\ref{PDPthBDCatAlg}; 
finally, the last equality follows from the fact that, according to Proposition~\ref{PDKerCHDCatHom}, $\mathrm{pr}_{\boldsymbol{\mathcal{B}}^{(2)}, \varphi}^{\llbracket\cdot\rrbracket}$ is a $\Sigma^{\boldsymbol{\mathcal{A}}^{(2)}}$-homomorphism from $\mathbf{Pth}_{\boldsymbol{\mathcal{B}}^{(2)}}^{\mathbf{f}^{(2)}}$ to $\llbracket\mathbf{Pth}_{\boldsymbol{\mathcal{B}}^{(2)}}^{\mathbf{f}^{(2)}}\rrbracket$.

Hence, $\mathrm{pr}_{\boldsymbol{\mathcal{B}}^{(2)}, \varphi}^{\llbracket\cdot\rrbracket} \circ f^{(2)\flat}$ is compatible with the $0$-source operation.

{\sffamily The mapping $\mathrm{pr}_{\boldsymbol{\mathcal{B}}^{(2)}, \varphi}^{\llbracket\cdot\rrbracket} \circ f^{(2)\flat}$ is compatible with the $0$-target.}

Let $s$ be a sort in $S$ and let us consider the $0$-target operation symbol $\mathrm{tg}_{s}^{0}$ in $\Sigma^{\boldsymbol{\mathcal{A}}^{(2)}}_{s,s}$. Let $\mathfrak{P}^{(2)}$ be a path in $\mathrm{Pth}_{\boldsymbol{\mathcal{B}}^{(2)}, \varphi(s)}$, then the following equality holds
$$
\left\llbracket
f_{s}^{(2)\flat} \left(
\mathrm{tg}_{s}^{0\mathbf{Pth}_{\boldsymbol{\mathcal{A}}^{(2)}}} \left(
\mathfrak{P}^{(2)}
\right)
\right)
\right\rrbracket_{\varphi(s)}
=
\mathrm{tg}_{s}^{0\llbracket\mathbf{Pth}_{\boldsymbol{\mathcal{B}}^{(2)}}^{\mathbf{f}^{(2)}}\rrbracket} \left(
\left\llbracket
f_{s}^{(2)\flat} \left(
\mathfrak{P}^{(2)}
\right)
\right\rrbracket_{\varphi(s)}
\right).
$$

The proof of this case is identical to that of the $0$-source.

Hence, $\mathrm{pr}_{\boldsymbol{\mathcal{B}}^{(2)}, \varphi}^{\llbracket\cdot\rrbracket} \circ f^{(2)\flat}$ is compatible with the $0$-target operation.

{\sffamily The mapping $\mathrm{pr}_{\boldsymbol{\mathcal{B}}^{(2)}, \varphi}^{\llbracket\cdot\rrbracket} \circ f^{(2)\flat}$ is compatible with the $0$-composition.}

Let $s$ be a sort in $S$ and let us consider the $0$-composition operation symbol $\circ_{s}^{0}$ in $\Sigma_{ss,s}^{\boldsymbol{\mathcal{A}}^{(2)}}$. Let $\mathfrak{P}^{(2)}$ and $\mathfrak{Q}^{(2)}$ be two paths in $\mathrm{Pth}_{\boldsymbol{\mathcal{A}}^{(2)}, s}$ such that
$$
\mathrm{sc}_{\boldsymbol{\mathcal{A}}^{(2)}, s}^{(0,2)}\left(\mathfrak{Q}^{(2)}\right)
=
\mathrm{tg}_{\boldsymbol{\mathcal{A}}^{(2)}, s}^{(0,2)}\left(\mathfrak{P}^{(2)}\right).
$$

We will prove that
$$
\left\llbracket
f_{s}^{(2)\flat}\left(
\mathfrak{Q}^{(2)}
\circ_{s}^{0\mathbf{Pth}_{\boldsymbol{\mathcal{A}}^{(2)}}}
\mathfrak{P}^{(2)}
\right)
\right\rrbracket_{\varphi(s)}
=
\left\llbracket
f_{s}^{(2)\flat}\left(
\mathfrak{Q}^{(2)}
\right)
\right\rrbracket_{\varphi(s)}
\circ_{s}^{0\llbracket\mathbf{Pth}_{\boldsymbol{\mathcal{B}}^{(2)}}^{\mathbf{f}^{(2)}}\rrbracket}
\left\llbracket
f_{s}^{(2)\flat}\left(
\mathfrak{P}^{(2)}
\right)
\right\rrbracket_{\varphi(s)}
$$

We will distinguish two cases according to the nature of the pair $\mathfrak{P}^{(2)}$ and $\mathfrak{Q}^{(2)}$. It could be the case that either $(1)$ both $\mathfrak{P}^{(2)}$ and $\mathfrak{Q}^{(2)}$ are $(2,[1])$-identity second-order paths or $(2)$ at least one of the second-order paths $\mathfrak{P}^{(2)}$ and $\mathfrak{Q}^{(2)}$ has length at least one.

If~$(1)$, i.e., if $\mathfrak{P}^{(2)}$ and $\mathfrak{Q}^{(2)}$ are $(2,[1])$-identity second-order paths, then $\mathfrak{P}^{(2)} = \mathrm{ip}_{\boldsymbol{\mathcal{A}}^{(2)}, s}^{(2,[1])\sharp} \left([P]_{s}\right)$ and $\mathfrak{Q}^{(2)} = \mathrm{ip}_{\boldsymbol{\mathcal{A}}^{(2)}, s}^{(2,[1])\sharp} \left([Q]_{s}\right)$ for some path term classes $[P]_{s}$ and $[Q]_{s}$ in $[\mathrm{PT}_{\boldsymbol{\mathcal{A}}^{(2)}}]_{s}$. Thus, according to Proposition~\ref{PDUIp}, it follows that
$$
\mathfrak{Q}^{(2)}
\circ_{s}^{0\mathbf{Pth}_{\boldsymbol{\mathcal{A}}^{(2)}}}
\mathfrak{P}^{(2)}
=
\mathrm{ip}_{\boldsymbol{\mathcal{A}}^{(2)}, s}^{(2,[1])\sharp}\left(
[Q \circ_{s}^{\mathbf{PT}_{\boldsymbol{\mathcal{A}}^{(1)}}} P]_{s}
\right).
$$

Therefore, the following chain of equalities holds
\begin{flushleft}
$
\left\llbracket
f_{s}^{(2)\flat}\left(
\mathfrak{Q}^{(2)}
\circ_{s}^{0\mathbf{Pth}_{\boldsymbol{\mathcal{A}}^{(2)}}}
\mathfrak{P}^{(2)}
\right)
\right\rrbracket_{\varphi(s)}
$
\allowdisplaybreaks
\begin{align}
&=
\left\llbracket
f_{s}^{(2)\flat}\left(
\mathrm{ip}_{\boldsymbol{\mathcal{A}}^{(2)}, s}^{(2,[1])\sharp}\left(
[Q \circ_{s}^{0\mathbf{PT}_{\boldsymbol{\mathcal{A}}^{(1)}}} P]_{s}
\right)
\right)
\right\rrbracket_{\varphi(s)}
\tag{1}
\\
&=
\left\llbracket
\mathrm{ip}_{\boldsymbol{\mathcal{B}}^{(2)}, \varphi(s)}^{(2,[1])\sharp}\left(
f_{s}^{(1)@}\left(
[Q \circ_{s}^{0\mathbf{PT}_{\boldsymbol{\mathcal{A}}^{(1)}}} P]_{s}
\right)
\right)
\right\rrbracket_{\varphi(s)}
\tag{2}
\\
&=
\left\llbracket
\mathrm{ip}_{\boldsymbol{\mathcal{B}}^{(2)}, \varphi(s)}^{(2,[1])\sharp}\left(
f_{s}^{(1)@}\left(
[Q]_{s} \circ_{s}^{0[\mathbf{PT}_{\boldsymbol{\mathcal{A}}^{(1)}}]} [P]_{s}
\right)
\right)
\right\rrbracket_{\varphi(s)}
\tag{3}
\\
&=
\left\llbracket
\mathrm{ip}_{\boldsymbol{\mathcal{B}}^{(2)}, \varphi(s)}^{(2,[1])\sharp}\left(
f_{s}^{(1)@}\left(
[Q]_{s}
\right)
\circ_{\varphi(s)}^{0[\mathbf{PT}_{\boldsymbol{\mathcal{B}}^{(1)}}]}
f_{s}^{(1)@}\left(
[P]_{s}
\right)
\right)
\right\rrbracket_{\varphi(s)}
\tag{4}
\\
&=
\left\llbracket
\mathrm{ip}_{\boldsymbol{\mathcal{B}}^{(2)}, \varphi(s)}^{(2,[1])\sharp}\left(
f_{s}^{(1)@}\left(
[Q]_{s}
\right)
\right)
\circ_{\varphi(s)}^{0\mathbf{Pth}_{\boldsymbol{\mathcal{B}}^{(2)}}}
\mathrm{ip}_{\boldsymbol{\mathcal{B}}^{(2)}, \varphi(s)}^{(2,[1])\sharp}\left(
f_{s}^{(1)@}\left(
[P]_{s}
\right)
\right)
\right\rrbracket_{\varphi(s)}
\tag{5}
\\
&=
\left\llbracket
f_{s}^{(2)\flat}\left(
\mathfrak{Q}^{(2)}
\right)
\circ_{\varphi(s)}^{0\mathbf{Pth}_{\boldsymbol{\mathcal{B}}^{(2)}}}
f_{s}^{(2)\flat}\left(
\mathfrak{P}^{(2)}
\right)
\right\rrbracket_{\varphi(s)}
\tag{6}
\\
&=
\left\llbracket
f_{s}^{(2)\flat}\left(
\mathfrak{Q}^{(2)}
\right)
\circ_{s}^{0\mathbf{Pth}_{\boldsymbol{\mathcal{B}}^{(2)}}^{\mathbf{f}^{(2)}}}
f_{s}^{(2)\flat}\left(
\mathfrak{P}^{(2)}
\right)
\right\rrbracket_{\varphi(s)}
\tag{7}
\\
&=
\left\llbracket
f_{s}^{(2)\flat}\left(
\mathfrak{Q}^{(2)}
\right)
\right\rrbracket_{\varphi(s)}
\circ_{s}^{0\llbracket\mathbf{Pth}_{\boldsymbol{\mathcal{B}}^{(2)}}^{\mathbf{f}^{(2)}}\rrbracket}
\left\llbracket
f_{s}^{(2)\flat}\left(
\mathfrak{P}^{(2)}
\right)
\right\rrbracket_{\varphi(s)}.
\tag{8}
\end{align}
\end{flushleft}

The first equality follows from Proposition~\ref{PDUIp};
the second equality unravels the definition of the mapping $f^{(2)\flat}$ at a $(2,[1])$-identity path;
the third equality follows from the definition of the partial $\Sigma^{\boldsymbol{\mathcal{A}}^{(1)}}$-algebra $[\mathbf{PT}_{\boldsymbol{\mathcal{A}}^{(1)}}]$;
the fourth equality follows from the fact that, according to Definition~\ref{DQPthExt}, $f^{[1]@}$ is a $\Sigma^{\boldsymbol{\mathcal{A}}^{(1)}}$-homomorphism;
the fifth equality follows from the fact that, according to Proposition~\ref{PDUIpCatHom}, $\mathrm{ip}^{(2,[1])\sharp}_{\boldsymbol{\mathcal{B}}^{(2)}, \varphi(s)}$ is a $\Sigma^{\boldsymbol{\mathcal{A}}^{(1)}}$-homomorphism;
the sixth equality recovers the definition of the mapping $f^{(2)\flat}$ at a $(2,[1])$-identity path;
the seventh equality follows from the definition of the operation symbol $\circ^{0}_{s}$ on the partial $\Sigma^{\boldsymbol{\mathcal{A}}^{(1)}}$-algebra $\mathbf{Pth}_{\boldsymbol{\mathcal{B}}^{(2)}}^{\mathbf{f}^{(2)}}$, introduced in Proposition~\ref{PDPthBCatAlg};
finally, the last equality follows from the definition of the operation symbol $\circ^{0}_{s}$ on the partial $\Sigma^{\boldsymbol{\mathcal{A}}^{(1)}}$-algebra $\llbracket\mathbf{Pth}_{\boldsymbol{\mathcal{B}}^{(2)}}^{\mathbf{f}^{(2)}}\rrbracket$, introduced in Proposition~\ref{PDQPthBDCatAlg}.

If~$(2)$, i.e., if at least one of the second-order paths $\mathfrak{P}^{(2)}$ or $\mathfrak{Q}^{(2)}$ has length at least one, then, by Corollary~\ref{CDPthWB}, we have that the second-order path $\mathfrak{Q}^{(2)}\circ_{s}^{0\mathbf{Pth}_{\boldsymbol{\mathcal{A}}^{(2)}}}\mathfrak{P}^{(2)}$ is a head-constant coherent and echelonless second-order path. Moreover, by Proposition~\ref{PDRecov}, the second-order path extraction procedure from Lemma~\ref{LDPthExtract} applied to the second-order path $\mathfrak{Q}^{(2)}\circ_{s}^{0\mathbf{Pth}_{\boldsymbol{\mathcal{A}}^{(2)}}}\mathfrak{P}^{(2)}$ retrieves the pair of second-order paths $(\mathfrak{Q}^{(2)}, \mathfrak{P}^{(2)})$.

Thus, according to Proposition~\ref{PDPthExt}, the value of the mapping $f^{(2)\flat}$ at $\mathfrak{Q}^{(2)}\circ_{s}^{0\mathbf{Pth}_{\boldsymbol{\mathcal{A}}^{(2)}}}\mathfrak{P}^{(2)}$ is given by
$$
f^{(2)\flat} \left(
\mathfrak{Q}^{(2)}
\circ_{s}^{0\mathbf{Pth}_{\boldsymbol{\mathcal{A}}^{(2)}}}
\mathfrak{P}^{(2)}
\right)
=
f^{(2)\flat} \left(
\mathfrak{Q}^{(2)}
\right)
\circ_{s}^{0\mathbf{Pth}_{\boldsymbol{\mathcal{B}}^{(2)}}^{\mathbf{f}^{(2)}}}
f^{(2)\flat} \left(
\mathfrak{P}^{(2)}
\right).
$$
Thus, the desired equality follows.

{\sffamily The mapping $\mathrm{pr}_{\boldsymbol{\mathcal{B}}^{(2)}, \varphi}^{\llbracket\cdot\rrbracket} \circ f^{(2)\flat}$ is compatible with the second-order rewrite rules.}

Let $s$ be a sort in $S$ and let $\mathfrak{p}^{(2)}$ be a second-order rewrite rule in $\mathcal{A}_{s}^{(2)}$. We have to show that
$$
\left\llbracket
f_{s}^{(2)\flat}\left(
\mathfrak{p}^{(2)\mathbf{Pth}_{\boldsymbol{\mathcal{A}}^{(2)}}}
\right)
\right\rrbracket_{\varphi(s)}
=
\mathfrak{p}^{(2)\llbracket\mathbf{Pth}_{\boldsymbol{\mathcal{B}}^{(2)}}^{\mathbf{f}^{(2)}}\rrbracket}
$$

The following chain of equalities holds
\allowdisplaybreaks
\begin{align*}
\left\llbracket
f_{s}^{(2)\flat}\left(
\mathfrak{p}^{(2)\mathbf{Pth}_{\boldsymbol{\mathcal{A}}^{(2)}}}
\right)
\right\rrbracket_{\varphi(s)}
&=
\left\llbracket
f_{s}^{(2)}\left(
\mathrm{ech}_{\boldsymbol{\mathcal{A}}^{(2)}, s}^{(2,\mathcal{A}^{(2)})}\left(
\mathfrak{p}^{(2)}
\right)
\right)
\right\rrbracket_{\varphi(s)}
\tag{1}
\\
&=
\left\llbracket
f_{s}^{(2)}\left(
\mathfrak{p}^{(2)}
\right)
\right\rrbracket_{\varphi(s)}
\tag{2}
\\
&=
\mathfrak{p}^{(2)\llbracket\mathbf{Pth}_{\boldsymbol{\mathcal{B}}^{(2)}}^{\mathbf{f}^{(2)}}\rrbracket}.
\tag{3}
\end{align*}

The first equality unravels the interpretation of the constant symbol $\mathfrak{p}$ in the partial $\Sigma^{\boldsymbol{\mathcal{A}}^{(2)}}$-algebra $\mathbf{Pth}_{\boldsymbol{\mathcal{A}}^{(2)}}$, introduced in Proposition~\ref{PDPthCatAlg};
the second equality unravels the definition of $f^{(2)\flat}$ at a second-order echelon path;
finally, the last equality recovers the interpretation of the constant symbol $\mathfrak{p}^{(2)}$ in the partial $\Sigma^{\boldsymbol{\mathcal{A}}^{(2)}}$-algebra $\llbracket\mathbf{Pth}_{\boldsymbol{\mathcal{B}}^{(2)}}^{\mathbf{f}^{(2)}}\rrbracket$, introduced in Proposition~\ref{PDQPthBDCatAlg}.

Hence, $\mathrm{pr}_{\boldsymbol{\mathcal{B}}^{(2)}, \varphi}^{\llbracket\cdot\rrbracket} \circ f^{(2)\flat}$ is compatible with the second-order rewrite rules.

{\sffamily The mapping $\mathrm{pr}_{\boldsymbol{\mathcal{B}}^{(2)}, \varphi}^{\llbracket\cdot\rrbracket} \circ f^{(2)\flat}$ is compatible with the $1$-source.}

Let $s$ be a sort in $S$ and let us consider the $1$-source operation symbol $\mathrm{sc}_{s}^{1}$ in $\Sigma^{\boldsymbol{\mathcal{A}}^{(2)}}_{s,s}$. Let $\mathfrak{P}^{(2)}$ be a second-order path in $\mathrm{Pth}_{\boldsymbol{\mathcal{A}}^{(2)}, s}$.

The following chain of equalities holds
\begin{flushleft}
$
\left\llbracket
f_{s}^{(2)\flat} \left(
\mathrm{sc}_{s}^{1\mathbf{Pth}_{\boldsymbol{\mathcal{A}}^{(2)}}} \left(
\mathfrak{P}^{(2)}
\right)
\right)
\right\rrbracket_{\varphi(s)}
$
\allowdisplaybreaks
\begin{align*}
&=
\left\llbracket
f_{s}^{(2)\flat} \left(
\mathrm{ip}_{\boldsymbol{\mathcal{A}}^{(2)}, s}^{(2,[1])\sharp} \left(
\mathrm{sc}_{\boldsymbol{\mathcal{A}}^{(2)}, s}^{([1],2)} \left(
\mathfrak{P}^{(2)}
\right)
\right)
\right)
\right\rrbracket_{\varphi(s)}
\tag{1}
\\
&=
\left\llbracket
\mathrm{ip}_{\boldsymbol{\mathcal{B}}^{(2)}, \varphi(s)}^{(2,[1])\sharp} \left(
f_{s}^{(1)\mathsf{q}} \left(
\mathrm{sc}_{\boldsymbol{\mathcal{A}}^{(2)}, s}^{([1],2)} \left(
\mathfrak{P}^{(2)}
\right)
\right)
\right)
\right\rrbracket_{\varphi(s)}
\tag{2}
\\
&=
\left\llbracket
\mathrm{ip}_{\boldsymbol{\mathcal{B}}^{(2)}, \varphi(s)}^{(2,[1])\sharp} \left(
\mathrm{sc}_{\boldsymbol{\mathcal{B}}^{(2)}, \varphi(s)}^{([1],2)} \left(
f_{s}^{(2)\flat} \left(
\mathfrak{P}^{(2)}
\right)
\right)
\right)
\right\rrbracket_{\varphi(s)}
\tag{3}
\\
&=
\left\llbracket
\mathrm{sc}_{\varphi(s)}^{1\mathbf{Pth}_{\boldsymbol{\mathcal{B}}^{(2)}}} \left(
f_{s}^{(2)\flat} \left(
\mathfrak{P}^{(2)}
\right)
\right)
\right\rrbracket_{\varphi(s)}
\tag{4}
\\
&=
\left\llbracket
\mathrm{sc}_{s}^{1\mathbf{Pth}_{\boldsymbol{\mathcal{B}}^{(2)}}^{\mathbf{f}^{(2)}}} \left(
f_{s}^{(2)\flat} \left(
\mathfrak{P}^{(2)}
\right)
\right)
\right\rrbracket_{\varphi(s)}
\tag{5}
\\
&=
\mathrm{sc}_{s}^{1\llbracket\mathbf{Pth}_{\boldsymbol{\mathcal{B}}^{(2)}}^{\mathbf{f}^{(2)}}\rrbracket} \left(
\left\llbracket
f_{s}^{(2)\flat} \left(
\mathfrak{P}^{(2)}
\right)
\right\rrbracket_{\varphi(s)}
\right).
\tag{6}
\end{align*}
\end{flushleft}

The first equality unravels the interpretation of the operation symbol $\mathrm{sc}_{s}^{1}$ in the partial $\Sigma^{\boldsymbol{\mathcal{A}}^{(2)}}$-algebra $\mathbf{Pth}_{\boldsymbol{\mathcal{A}}^{(2)}}$, introduced in Proposition~\ref{PDPthDCatAlg}; 
the second equality follows from Proposition~\ref{PDPthExt}; 
the third equality follows from Proposition~\ref{PDPthExt}; 
the fourth equality recovers the interpretation of the operation symbol $\mathrm{sc}_{\varphi(s)}^{1}$ in the partial $\Lambda^{\boldsymbol{\mathcal{B}}^{(2)}}$-algebra $\mathbf{Pth}_{\boldsymbol{\mathcal{B}}^{(2)}}$, according to Proposition~\ref{PDPthDCatAlg}; 
the fifth equality recovers the interpretation of the operation symbol $\mathrm{sc}_{s}^{1}$ in the partial $\Sigma^{\boldsymbol{\mathcal{A}}^{(2)}}$-algebra $\mathbf{Pth}_{\boldsymbol{\mathcal{B}}^{(2)}}^{\mathbf{f}^{(2)}}$, introduced in Proposition~\ref{PDPthBDCatAlg}; 
finally, the last equality follows from the fact that, according to Proposition~\ref{PDKerCHDCatHom}, $\mathrm{pr}_{\boldsymbol{\mathcal{B}}^{(2)}, \varphi}^{\llbracket\cdot\rrbracket}$ is a $\Sigma^{\boldsymbol{\mathcal{A}}^{(2)}}$-homomorphism from $\mathbf{Pth}_{\boldsymbol{\mathcal{B}}^{(2)}}^{\mathbf{f}^{(2)}}$ to $\llbracket\mathbf{Pth}_{\boldsymbol{\mathcal{B}}^{(2)}}^{\mathbf{f}^{(2)}}\rrbracket$.

Hence, $\mathrm{pr}_{\boldsymbol{\mathcal{B}}^{(2)}, \varphi}^{\llbracket\cdot\rrbracket} \circ f^{(2)\flat}$ is compatible with the $1$-source operation.

{\sffamily The mapping $\mathrm{pr}_{\boldsymbol{\mathcal{B}}^{(2)}, \varphi}^{\llbracket\cdot\rrbracket} \circ f^{(2)\flat}$ is compatible with the $1$-target.}

Let $s$ be a sort in $S$ and let us consider the $1$-target operation symbol $\mathrm{tg}_{s}^{1}$ in $\Sigma^{\boldsymbol{\mathcal{A}}^{(2)}}_{s,s}$. Let $\mathfrak{P}^{(2)}$ be a path in $\mathrm{Pth}_{\boldsymbol{\mathcal{B}}^{(2)}, \varphi(s)}$, then the following equality holds
$$
\left\llbracket
f_{s}^{(2)\flat} \left(
\mathrm{tg}_{s}^{1\mathbf{Pth}_{\boldsymbol{\mathcal{A}}^{(2)}}} \left(
\mathfrak{P}^{(2)}
\right)
\right)
\right\rrbracket_{\varphi(s)}
=
\mathrm{tg}_{s}^{1\llbracket\mathbf{Pth}_{\boldsymbol{\mathcal{B}}^{(2)}}^{\mathbf{f}^{(2)}}\rrbracket} \left(
\left\llbracket
f_{s}^{(2)\flat} \left(
\mathfrak{P}^{(2)}
\right)
\right\rrbracket_{\varphi(s)}
\right).
$$

The proof of this case is identical to that of the $1$-source.

Hence, $\mathrm{pr}_{\boldsymbol{\mathcal{B}}^{(2)}, \varphi}^{\llbracket\cdot\rrbracket} \circ f^{(2)\flat}$ is compatible with the $1$-target operation.

{\sffamily The mapping $\mathrm{pr}_{\boldsymbol{\mathcal{B}}^{(2)}, \varphi}^{\llbracket\cdot\rrbracket} \circ f^{(2)\flat}$ is compatible with the $1$-composition.}

Let $s$ be a sort in $S$ and let us consider the $1$-composition operation symbol $\circ_{s}^{1}$ in $\Sigma_{ss,s}^{\boldsymbol{\mathcal{A}}^{(2)}}$. Let $\mathfrak{P}^{(2)}$ and $\mathfrak{Q}^{(2)}$ be two paths in $\mathrm{Pth}_{\boldsymbol{\mathcal{A}}^{(2)}, s}$ such that
$$
\mathrm{sc}_{\boldsymbol{\mathcal{A}}^{(2)}, s}^{([1],2)}\left(\mathfrak{Q}^{(2)}\right)
=
\mathrm{tg}_{\boldsymbol{\mathcal{A}}^{(2)}, s}^{([1],2)}\left(\mathfrak{P}^{(2)}\right).
$$

We will prove that
$$
\left\llbracket
f_{s}^{(2)\flat}\left(
\mathfrak{Q}^{(2)}
\circ_{s}^{1\mathbf{Pth}_{\boldsymbol{\mathcal{A}}^{(2)}}}
\mathfrak{P}^{(2)}
\right)
\right\rrbracket_{\varphi(s)}
=
\left\llbracket
f_{s}^{(2)\flat}\left(
\mathfrak{Q}^{(2)}
\right)
\right\rrbracket_{\varphi(s)}
\circ_{s}^{1\llbracket\mathbf{Pth}_{\boldsymbol{\mathcal{B}}^{(2)}}^{\mathbf{f}^{(2)}}\rrbracket}
\left\llbracket
f_{s}^{(2)\flat}\left(
\mathfrak{P}^{(2)}
\right)
\right\rrbracket_{\varphi(s)}
$$
by Artinian induction of $(\coprod \mathrm{Pth}_{\boldsymbol{\mathcal{A}}^{(2)}}, \leq_{\mathrm{Pth}_{\boldsymbol{\mathcal{A}}^{(2)}}})$.

{\sf Base step of the Artinian induction.}

Let $(\mathfrak{Q}^{(2)}\circ_{s}^{1\mathbf{Pth}_{\boldsymbol{\mathcal{A}}^{(2)}}}\mathfrak{P}^{(2)}, s)$ be a minimal element in $(\coprod \mathrm{Pth}_{\boldsymbol{\mathcal{A}}^{(2)}}, \leq_{\mathrm{Pth}_{\boldsymbol{\mathcal{A}}^{(2)}}})$. Then, by Proposition~\ref{PDMinimal}, the second order path $\mathfrak{Q}^{(2)}\circ_{s}^{1\mathbf{Pth}_{\boldsymbol{\mathcal{A}}^{(2)}}}\mathfrak{P}^{(2)}$ is either a $(2,[1])$-identity second-order path or a second-order echelon. In any case, either $\mathfrak{Q}^{(2)}$ or $\mathfrak{P}^{(2)}$ must be a $(2,[1])$-identity second-order path. 

Assume that $\mathfrak{P}^{(2)}$ is a $(2,[1])$-identity path.

Since $\mathfrak{Q}^{(2)}\circ_{s}^{1\mathbf{Pth}_{\boldsymbol{\mathcal{A}}^{(2)}}}\mathfrak{P}^{(2)}$ is defined, we have that $\mathrm{sc}_{\boldsymbol{\mathcal{A}}^{(2)}, s}^{([1],2)}(\mathfrak{Q}^{(2)}) = \mathrm{tg}_{\boldsymbol{\mathcal{A}}^{(2)}, s}^{([1],2)}(\mathfrak{P}^{(2)})$. Hence, $\mathfrak{P}^{(2)}$ is the $(2,[1])$-identity second-order path on the $([1],2)$-source of $\mathfrak{Q}^{(2)}$, i.e., $\mathfrak{P}^{(2)} = \mathrm{ip}_{\boldsymbol{\mathcal{A}}^{(2)}, s}^{(2,[1])\sharp}(\mathrm{sc}_{s}^{([1],2)}(\mathfrak{P}^{(2)}))$. Hence, the $1$-composition $\mathfrak{Q}^{(2)} \circ_{s}^{1\mathbf{Pth}_{\boldsymbol{\mathcal{A}}^{(2)}}}\mathfrak{P}^{(2)}$ reduces to $\mathfrak{Q}^{(2)}$.

On the other hand, the following chain of equalities holds
\begin{align}
\left\llbracket
f_{s}^{(2)\flat}\left(
\mathfrak{P}^{(2)}
\right)
\right\rrbracket_{\varphi(s)}
&=
\left\llbracket
f_{s}^{(2)\flat}\left(
\mathrm{ip}^{(2,[1])\sharp}_{\boldsymbol{\mathcal{A}}^{(2)},s}\left(
\mathrm{sc}^{([1],2)}_{\boldsymbol{\mathcal{A}}^{(2)},s}\left(
\mathfrak{Q}^{(2)}
\right)
\right)
\right)
\right\rrbracket_{\varphi(s)}
\tag{1}
\\
&=
\left\llbracket
\mathrm{ip}^{(2,[1])\sharp}_{\boldsymbol{\mathcal{B}}^{(2)},\varphi(s)}\left(
f_{s}^{(1)@}\left(
\mathrm{sc}^{([1],2)}_{\boldsymbol{\mathcal{A}}^{(2)},s}\left(
\mathfrak{Q}^{(2)}
\right)
\right)
\right)
\right\rrbracket_{\varphi(s)}
\tag{2}
\\
&=
\left\llbracket
\mathrm{ip}^{(2,[1])\sharp}_{\boldsymbol{\mathcal{B}}^{(2)},\varphi(s)}\left(
\mathrm{sc}^{([1],2)}_{\boldsymbol{\mathcal{B}}^{(2)},\varphi(s)}\left(
f_{s}^{(2)\sharp}\left(
\mathfrak{Q}^{(2)}
\right)
\right)
\right)
\right\rrbracket_{\varphi(s)}.
\tag{3}
\end{align}

The first equality follows from the fact that $\mathfrak{P}^{(2)} = \mathrm{ip}^{(2,[1])\sharp}_{\boldsymbol{\mathcal{A}}^{(2)}, s}(\mathrm{sc}^{([1],2)}_{\boldsymbol{\mathcal{A}}^{(2)}, s}(\mathfrak{Q}^{(2)}))$;
the second-equality follows from Proposition~\ref{PDPthExt};
finally, the last equality follows from Proposition~\ref{PDPthExt}.

Hence, the $1$-composition 
$
\llbracket
f_{s}^{(2)\flat}(
\mathfrak{Q}^{(2)}
)
\rrbracket_{\varphi(s)}
\circ_{s}^{1\llbracket\mathbf{Pth}_{\boldsymbol{\mathcal{B}}^{(2)}}^{\mathbf{f}^{(2)}}\rrbracket}
\llbracket
f_{s}^{(2)\flat}(
\mathfrak{P}^{(2)}
)
\rrbracket_{\varphi(s)}
$
reduces to
$
\llbracket
f_{s}^{(2)\flat}(
\mathfrak{Q}^{(2)}
)
\rrbracket_{\varphi(s)}
$.

All in all, we conclude that
\begin{flushleft}
$
\left\llbracket
f_{s}^{(2)\flat}\left(
\mathfrak{Q}^{(2)}
\circ_{s}^{1\mathbf{Pth}_{\boldsymbol{\mathcal{A}}^{(2)}}}
\mathfrak{P}^{(2)}
\right)
\right\rrbracket_{\varphi(s)}
$
\allowdisplaybreaks
\begin{align}
&=
\left\llbracket
f_{s}^{(2)\flat}\left(
\mathfrak{Q}^{(2)}
\right)
\right\rrbracket_{\varphi(s)}
\tag{1}
\\
&=
\left\llbracket
f_{s}^{(2)\flat}\left(
\mathfrak{Q}^{(2)}
\right)
\right\rrbracket_{\varphi(s)}
\circ_{s}^{1\llbracket\mathbf{Pth}_{\boldsymbol{\mathcal{B}}^{(2)}}^{\mathbf{f}^{(2)}}\rrbracket}
\left\llbracket
f_{s}^{(2)\flat}\left(
\mathfrak{P}^{(2)}
\right)
\right\rrbracket_{\varphi(s)}
\tag{2}
\end{align}
\end{flushleft}

The first equality follows from the fact that the $1$-composition $\mathfrak{Q}^{(2)} \circ_{s}^{1\mathbf{Pth}_{\boldsymbol{\mathcal{A}}^{(2)}}}\mathfrak{P}^{(2)}$ reduces to $\mathfrak{Q}^{(2)}$; 
the second equality follows from the fac that the $1$-composition 
$
\llbracket
f_{s}^{(2)\flat}(
\mathfrak{Q}^{(2)}
)
\rrbracket_{\varphi(s)}
\circ_{s}^{1\llbracket\mathbf{Pth}_{\boldsymbol{\mathcal{B}}^{(2)}}^{\mathbf{f}^{(2)}}\rrbracket}
\llbracket
f_{s}^{(2)\flat}(
\mathfrak{P}^{(2)}
)
\rrbracket_{\varphi(s)}
$
reduces to
$
\llbracket
f_{s}^{(2)\flat}(
\mathfrak{Q}^{(2)}
)
\rrbracket_{\varphi(s)}
$.

The case where $\mathfrak{Q}^{(2)}$ is a $(2,[1])$-identity path follows similarly.

This concludes the base step of the Artinian induction.

{\sf Inductive step of the Artinian induction.}

Let $(\mathfrak{Q}^{(2)} \circ_{s}^{1\mathbf{Pth}_{\boldsymbol{\mathcal{A}}^{(2)}}} \mathfrak{P}^{(2)}, s)$ be a non-minimal element in $(\coprod\mathrm{Pth}_{\boldsymbol{\mathcal{A}}^{(2)}}, \leq_{\mathbf{Pth}_{\boldsymbol{\mathcal{A}}^{(2)}}})$. Let us suppose that, for every sort $t \in S$ and every path $\mathfrak{Q}'^{(2)} \circ_{t}^{1\mathbf{Pth}_{\boldsymbol{\mathcal{A}}^{(2)}}} \mathfrak{P}'^{(2)} \in \mathrm{Pth}_{\boldsymbol{\mathcal{A}}^{(2)}, t}$, if $(\mathfrak{Q}'^{(2)} \circ_{t}^{1\mathbf{Pth}_{\boldsymbol{\mathcal{A}}^{(2)}}} \mathfrak{P}'^{(2)}, t) <_{\mathbf{Pth}_{\boldsymbol{\mathcal{A}}^{(2)}}} (\mathfrak{Q}^{(2)} \circ_{s}^{1\mathbf{Pth}_{\boldsymbol{\mathcal{A}}^{(2)}}} \mathfrak{P}^{(2)}, s)$, then
$$
\left\llbracket
f_{t}^{(2)\flat}\left(
\mathfrak{Q}'^{(2)}
\circ_{t}^{1\mathbf{Pth}_{\boldsymbol{\mathcal{A}}^{(2)}}}
\mathfrak{P}'^{(2)}
\right)
\right\rrbracket_{\varphi(t)}
=
\left\llbracket
f_{t}^{(2)\flat}\left(
\mathfrak{Q}'^{(2)}
\right)
\right\rrbracket_{\varphi(t)}
\circ_{t}^{1\llbracket\mathbf{Pth}_{\boldsymbol{\mathcal{B}}^{(2)}}^{\mathbf{f}^{(2)}}\rrbracket} 
\left\llbracket
f_{t}^{(2)\flat}\left(
\mathfrak{P}'^{(2)}
\right)
\right\rrbracket_{\varphi(t)}.
$$

Since $\mathfrak{Q}^{(2)} \circ_{s}^{1\mathbf{Pth}_{\boldsymbol{\mathcal{A}}^{(2)}}} \mathfrak{P}^{(2)}$ is a non-minimal element in $(\coprod\mathrm{Pth}_{\boldsymbol{\mathcal{A}}^{(2)}}, \leq_{\mathbf{Pth}_{\boldsymbol{\mathcal{A}}^{(2)}}})$ and we can assume that neither $\mathfrak{Q}^{(2)}$ nor $\mathfrak{P}^{(2)}$ are $(2,[1])$-identity paths because this case already has been considered, we have, by Lemma~\ref{LDOrdI}, that $\mathfrak{Q}^{(2)} \circ_{s}^{1\mathbf{Pth}_{\boldsymbol{\mathcal{A}}^{(2)}}} \mathfrak{P}^{(2)}$ is either~(1) a second-order path of length strictly grater than one containing at least one second-order echelon or~(2) an echelonless second-order path.

If~(1), then let $i\in\bb{\mathfrak{Q}^{(2)} \circ_{s}^{1\mathbf{Pth}_{\boldsymbol{\mathcal{A}}^{(2)}}} \mathfrak{P}^{(2)}}$ be the first index for which the one-step subpath $(\mathfrak{Q}^{(2)} \circ_{s}^{1\mathbf{Pth}_{\boldsymbol{\mathcal{A}}^{(2)}}} \mathfrak{P}^{(2)})^{i,i}$ of $\mathfrak{Q}^{(2)} \circ_{s}^{1\mathbf{Pth}_{\boldsymbol{\mathcal{A}}^{(2)}}} \mathfrak{P}^{(2)}$ is a second-order echelon. We distinguish the cases~(1.1) $i=0$ and~(1.2) $i>0$.

If~(1.1), i.e., $i=0$, since we are assuming that $\mathfrak{P}^{(2)}$ is not a (2,[1])-identity second-order path, we have that $\mathfrak{P}^{(2)}$ has a second-order echelon on its first step. Then it could be the case that either~(1.1.1) $\mathfrak{P}^{(2)}$ is a second-order echelon or~(1.1.2) $\mathfrak{P}^{(2)}$ is a second-order path of length strictly greater than one containing a second-order echelon on its first step.

If~(1.1.1), then, regarding the second-order paths $\mathfrak{Q}^{(2)}$ and $\mathfrak{P}^{(2)}$, we have that
\begin{enumerate}
\item[(i)]
$\mathfrak{P}^{(2)}$ is a second-order echelon.
\item[(ii)]
$\mathfrak{Q}^{(2)}$ is a second-order path.
\end{enumerate}

So, considering the foregoing, we can affirm that
\begin{flushleft}
$
\left\llbracket
f_{s}^{(2)\flat}\left(
\mathfrak{Q}^{(2)}
\circ_{s}^{1\mathbf{Pth}_{\boldsymbol{\mathcal{A}}^{(2)}}}
\mathfrak{P}^{(2)}
\right)
\right\rrbracket_{\varphi(s)}
$
\allowdisplaybreaks
\begin{align*}
&=
\left\llbracket
f_{s}^{(2)\flat}\left(
\mathfrak{Q}^{(2)}
\right)
\circ_{\varphi(s)}^{1\mathbf{Pth}_{\boldsymbol{\mathcal{B}}^{(2)}}}
f_{s}^{(2)\flat}\left(
\mathfrak{P}^{(2)}
\right)
\right\rrbracket_{\varphi(s)}
\tag{1}
\\
&=
\left\llbracket
f_{s}^{(2)\flat}\left(
\mathfrak{Q}^{(2)}
\right)
\circ_{s}^{1\mathbf{Pth}_{\boldsymbol{\mathcal{B}}^{(2)}}^{\mathbf{f}^{(2)}}}
f_{s}^{(2)\flat}\left(
\mathfrak{P}^{(2)}
\right)
\right\rrbracket_{\varphi(s)}
\tag{2}
\\
&=
\left\llbracket
f_{s}^{(2)\flat}\left(
\mathfrak{Q}^{(2)}
\right)
\right\rrbracket_{\varphi(s)}
\circ_{s}^{1\llbracket\mathbf{Pth}_{\boldsymbol{\mathcal{B}}^{(2)}}^{\mathbf{f}^{(2)}}\rrbracket}
\left\llbracket
f_{s}^{(2)\flat}\left(
\mathfrak{P}^{(2)}
\right)
\right\rrbracket_{\varphi(s)}.
\tag{3}
\end{align*}
\end{flushleft}

The first equality unravels the definition of the path extension mapping $f^{(2)\flat}$ introduced in Proposition~\ref{PDPthExt};
the second equality recovers the interpretation of the operation symbol $\circ_{s}^{1}$ in the partial $\Sigma^{\boldsymbol{\mathcal{A}}^{(2)}}$-algebra $\mathbf{Pth}_{\boldsymbol{\mathcal{B}}^{(2)}}^{\mathbf{f}^{(2)}}$ introduced in Proposition~\ref{PDPthBDCatAlg};
finally, the last equality recovers the interpretation of the operation symbol $\circ_{s}^{1}$ in the partial $\Sigma^{\boldsymbol{\mathcal{A}}^{(2)}}$-algebra $\llbracket\mathbf{Pth}_{\boldsymbol{\mathcal{B}}^{(2)}}^{\mathbf{f}^{(2)}}\rrbracket$ introduced in Proposition~\ref{PDQPthBDCatAlg}.

If~(1.1.2), then, regarding the second-order paths $\mathfrak{Q}^{(2)}$ and $\mathfrak{P}^{(2)}$, we have that
\begin{enumerate}
\item[(i)]
$\mathfrak{P}^{(2)}$ is a second-order path of length strictly grater than one containing a second-order echelon on its first step.
\item[(ii)]
$\mathfrak{Q}^{(2)}$ is a second-order path.
\end{enumerate}

From (i) and taking into account the definition of $f^{(2)\flat}$ introduced in Proposition~\ref{PDPthExt}, we have that the value of the second-order path extension mapping $f^{(2)\flat}$ at $\mathfrak{P}^{(2)}$ is given by
$$
f_{s}^{(2)\flat}\left(
\mathfrak{P}^{(2)}
\right)
=
f_{s}^{(2)\flat}\left(
\mathfrak{P}^{(2)1,\bb{\mathfrak{P}^{(2)}}-1}
\right)
\circ^{1\mathbf{Pth}_{\boldsymbol{\mathcal{B}}^{(2)}}}_{\varphi(s)}
f_{s}^{(2)\flat}\left(
\mathfrak{P}^{(2)0,0}
\right).
$$

Since $(\mathfrak{Q}^{(2)} \circ_{s}^{1\mathbf{Pth}_{\boldsymbol{\mathcal{A}}^{(2)}}} \mathfrak{P}^{(2)1,\bb{\mathfrak{P}^{(2)}}-1}, s)$$\prec_{\mathbf{Pth}_{\boldsymbol{\mathcal{A}}^{(2)}}}$$(\mathfrak{Q}^{(2)} \circ_{s}^{1\mathbf{Pth}_{\boldsymbol{\mathcal{A}}^{(2)}}} \mathfrak{P}^{(2)}, s)$ we have, by induction, that
\begin{align*}
&\left\llbracket
f_{s}^{(2)\flat}\left(
\mathfrak{Q}^{(2)}
\circ_{s}^{1\mathbf{Pth}_{\boldsymbol{\mathcal{A}}^{(2)}}}
\mathfrak{P}^{(2)1,\bb{\mathfrak{P}^{(2)}}-1}
\right)
\right\rrbracket_{\varphi(s)}
\\
&\hspace{3cm}
=
\left\llbracket
f_{s}^{(2)\flat}\left(
\mathfrak{Q}^{(2)}
\right)
\right\rrbracket_{\varphi(s)}
\circ_{s}^{1\llbracket\mathbf{Pth}_{\boldsymbol{\mathcal{B}}^{(2)}}^{\mathbf{f}^{(2)}}\rrbracket}
\left\llbracket
f_{s}^{(2)\flat}\left(
\mathfrak{P}^{(2)1,\bb{\mathfrak{P}^{(2)}}-1}
\right)
\right\rrbracket_{\varphi(s)}.
\end{align*}

So, considering the foregoing, we can affirm that
\begin{flushleft}
$
\left\llbracket
f_{s}^{(2)\flat}\left(
\mathfrak{Q}^{(2)} \circ_{s}^{1\mathbf{Pth}_{\boldsymbol{\mathcal{A}}^{(2)}}} \mathfrak{P}^{(2)}
\right)
\right\rrbracket_{\varphi(s)}
$
\allowdisplaybreaks
\begin{align*}
&=
\left\llbracket
f_{s}^{(2)\flat}\left(
\mathfrak{Q}^{(2)} \circ_{s}^{1\mathbf{Pth}_{\boldsymbol{\mathcal{A}}^{(2)}}} \mathfrak{P}^{(2)1,\bb{\mathfrak{P}^{(2)}}-1}
\right)
\circ_{\varphi(s)}^{1\mathbf{Pth}_{\boldsymbol{\mathcal{B}}^{(2)}}}
f_{s}^{(2)\flat}\left(
\mathfrak{P}^{(2)0,0}
\right)
\right\rrbracket_{\varphi(s)}
\tag{1}
\\
&=
\left\llbracket
f_{s}^{(2)\flat}\left(
\mathfrak{Q}^{(2)} \circ_{s}^{1\mathbf{Pth}_{\boldsymbol{\mathcal{A}}^{(2)}}} \mathfrak{P}^{(2)1,\bb{\mathfrak{P}^{(2)}}-1}
\right)
\circ_{s}^{1\mathbf{Pth}_{\boldsymbol{\mathcal{B}}^{(2)}}^{\mathbf{f}^{(2)}}}
f_{s}^{(2)\flat}\left(
\mathfrak{P}^{(2)0,0}
\right)
\right\rrbracket_{\varphi(s)}
\tag{2}
\\
&=
\left\llbracket
f_{s}^{(2)\flat}\left(
\mathfrak{Q}^{(2)} \circ_{s}^{1\mathbf{Pth}_{\boldsymbol{\mathcal{A}}^{(2)}}} \mathfrak{P}^{(2)1,\bb{\mathfrak{P}^{(2)}}-1}
\right)
\right\rrbracket_{\varphi(s)}
\circ_{s}^{1\llbracket\mathbf{Pth}_{\boldsymbol{\mathcal{B}}^{(2)}}^{\mathbf{f}^{(2)}}\rrbracket}
\left\llbracket
f_{s}^{(2)\flat}\left(
\mathfrak{P}^{(2)0,0}
\right)
\right\rrbracket_{\varphi(s)}
\tag{3}
\\
&=
\left(
\left\llbracket
f_{s}^{(2)\flat}\left(
\mathfrak{Q}^{(2)}
\right)
\right\rrbracket_{\varphi(s)}
\circ_{s}^{1\llbracket\mathbf{Pth}_{\boldsymbol{\mathcal{B}}^{(2)}}^{\mathbf{f}^{(2)}}\rrbracket} 
\left\llbracket
f_{s}^{(2)\flat}\left(
\mathfrak{P}^{(2)1,\bb{\mathfrak{P}}-1}
\right)
\right\rrbracket_{\varphi(s)}
\right)
\\
&\hspace{6.7cm}
\circ_{s}^{1\llbracket\mathbf{Pth}_{\boldsymbol{\mathcal{B}}^{(2)}}^{\mathbf{f}^{(2)}}\rrbracket}
\left\llbracket
f_{s}^{(2)\flat}\left(
\mathfrak{P}^{(2)0,0}
\right)
\right\rrbracket_{\varphi(s)}
\tag{4}
\\
&=
\left\llbracket
f_{s}^{(2)\flat}\left(
\mathfrak{Q}^{(2)}
\right)
\right\rrbracket_{\varphi(s)}
\circ_{s}^{1\llbracket\mathbf{Pth}_{\boldsymbol{\mathcal{B}}^{(2)}}^{\mathbf{f}^{(2)}}\rrbracket} 
\\
&\hspace{2.4cm}
\left(
\left\llbracket
f_{s}^{(2)\flat}\left(
\mathfrak{P}^{(2)1,\bb{\mathfrak{P}}-1}
\right)
\right\rrbracket_{\varphi(s)}
\circ_{s}^{1\llbracket\mathbf{Pth}_{\boldsymbol{\mathcal{B}}^{(2)}}^{\mathbf{f}^{(2)}}\rrbracket}
\left\llbracket
f_{s}^{(2)\flat}\left(
\mathfrak{P}^{(2)0,0}
\right)
\right\rrbracket_{\varphi(s)}
\right)
\tag{5}
\\
&=
\left\llbracket
f_{s}^{(2)\flat}\left(
\mathfrak{Q}^{(2)}
\right)
\right\rrbracket_{\varphi(s)}
\circ_{s}^{1\llbracket\mathbf{Pth}_{\boldsymbol{\mathcal{B}}^{(2)}}^{\mathbf{f}^{(2)}}\rrbracket}
\\
&\hspace{4cm}
\left\llbracket
f_{s}^{(2)\flat}\left(
\mathfrak{P}^{(2)1,\bb{\mathfrak{P}}-1}
\right)
\circ_{s}^{1\mathbf{Pth}_{\boldsymbol{\mathcal{B}}^{(2)}}^{\mathbf{f}^{(2)}}}
f_{s}^{(2)\flat}\left(
\mathfrak{P}^{(2)0,0}
\right)
\right\rrbracket_{\varphi(s)}
\tag{6}
\\
&=
\left\llbracket
f_{s}^{(2)\flat}\left(
\mathfrak{Q}^{(2)}
\right)
\right\rrbracket_{\varphi(s)}
\circ_{s}^{1\llbracket\mathbf{Pth}_{\boldsymbol{\mathcal{B}}^{(2)}}^{\mathbf{f}^{(2)}}\rrbracket} 
\\
&\hspace{4.1cm}
\left\llbracket
f_{s}^{(2)\flat}\left(
\mathfrak{P}^{(2)1,\bb{\mathfrak{P}}-1}
\right)
\circ_{\varphi(s)}^{1\mathbf{Pth}_{\boldsymbol{\mathcal{B}}^{(2)}}}
f_{s}^{(2)\flat}\left(
\mathfrak{P}^{(2)0,0}
\right)
\right\rrbracket_{\varphi(s)}
\tag{7}
\\
&=
\left\llbracket
f_{s}^{(2)\flat}\left(
\mathfrak{Q}^{(2)}
\right)
\right\rrbracket_{\varphi(s)}
\circ_{s}^{1\llbracket\mathbf{Pth}_{\boldsymbol{\mathcal{B}}^{(2)}}^{\mathbf{f}^{(2)}}\rrbracket} 
\left\llbracket
f_{s}^{(2)\flat}\left(
\mathfrak{P}^{(2)1,\bb{\mathfrak{P}}-1}
\circ_{s}^{0\mathbf{Pth}_{\boldsymbol{\mathcal{A}}^{(1)}}}
\mathfrak{P}^{(2)0,0}
\right)
\right\rrbracket_{\varphi(s)}
\tag{8}
\\
&=
\left\llbracket
f_{s}^{(2)\flat}\left(
\mathfrak{Q}^{(2)}
\right)
\right\rrbracket_{\varphi(s)}
\circ_{s}^{1\llbracket\mathbf{Pth}_{\boldsymbol{\mathcal{B}}^{(2)}}^{\mathbf{f}^{(2)}}\rrbracket} 
\left\llbracket
f_{s}^{(2)\flat}\left(
\mathfrak{P}^{(2)}
\right)
\right\rrbracket_{\varphi(s)}.
\tag{9}
\end{align*}
\end{flushleft}

The first equality unravels the definition of the second-order path extension mapping $f^{(2)\flat}$ introduced in Proposition~\ref{PDPthExt};
the second equality recovers the definition of the operation symbol $\circ_{s}^{1}$ in the partial $\Sigma^{\boldsymbol{\mathcal{A}}^{(2)}}$-algebra $\mathbf{Pth}_{\boldsymbol{\mathcal{B}}^{(2)}}^{\mathbf{f}^{(2)}}$ introduced in Proposition~\ref{PDPthBCatAlg};
the third equality recovers the definition of the operation symbol $\circ_{s}^{1}$ in the partial $\Sigma^{\boldsymbol{\mathcal{A}}^{(1)}}$-algebra $\llbracket\mathbf{Pth}_{\boldsymbol{\mathcal{B}}^{(2)}}^{\mathbf{f}^{(2)}}\rrbracket$ introduced in Proposition~\ref{PDQPthBCatAlg};
the fourth equality follows by Artinian induction;
the fifth equality follows from the fact that, according to Proposition~\ref{PDPthComp}, the interpretation of the partial operation $\circ_{s}^{1}$ is associative;
the sixth equality unravels the definition of the operation symbol $\circ_{s}^{1}$ in the partial $\Sigma^{\boldsymbol{\mathcal{A}}^{(2)}}$-algebra $\llbracket\mathbf{Pth}_{\boldsymbol{\mathcal{B}}^{(2)}}^{\mathbf{f}^{(2)}}\rrbracket$ introduced in Proposition~\ref{PDQPthBCatAlg};
The seventh equality unravels the definition of the operation symbol $\circ_{s}^{1}$ in the partial $\Sigma^{\boldsymbol{\mathcal{A}}^{(2)}}$-algebra $\mathbf{Pth}_{\boldsymbol{\mathcal{B}}^{(2)}}^{\mathbf{f}^{(2)}}$ introduced in Proposition~\ref{PDPthBCatAlg};
the eight equality recovers the definition of $f^{(2)\flat}$ since $\mathfrak{P}^{(2)}$ is a second-order path of length strictly greater than one containing a second-order echelon on its first step;
finally, the last equality recovers the definition of $\mathfrak{P}^{(2)}$.

This finishes the case $i = 0$.

For the case~$(1.2)$, i.e., if $i\neq 0$, since $\bb{\mathfrak{Q}^{(2)} \circ_{s}^{1\mathbf{Pth}_{\boldsymbol{\mathcal{A}}^{(2)}}} \mathfrak{P}^{(2)}} = \bb{\mathfrak{Q}^{(2)}} + \bb{\mathfrak{P}^{(2)}}$, then either~$(1.2.1)$ $i\in\bb{\mathfrak{P}^{(2)}}$ or~$(1.2.2)$ $i\in[\bb{P^{(2)}}, \bb{\mathfrak{Q}^{(2)} \circ_{s}^{1\mathbf{Pth}_{\boldsymbol{\mathcal{A}}^{(2)}}} \mathfrak{P}^{(2)}}-1]$.

If~(1.2.1), i.e., if we find ourselves in the case where $\mathfrak{Q}^{(2)} \circ_{s}^{1\mathbf{Pth}_{\boldsymbol{\mathcal{A}}^{(2)}}} \mathfrak{P}^{(2)}$ is a second-order path that is not echelonless and $i > 0$, $i \in \bb{\mathfrak{P}^{(2)}}$ is the first index for which the one step subpath $(\mathfrak{Q}^{(2)} \circ_{s}^{1\mathbf{Pth}_{\boldsymbol{\mathcal{A}}^{(2)}}} \mathfrak{P}^{(2)})^{i,i}$ is a second-order echelon then, regarding the second-order paths $\mathfrak{Q}^{(2)}$ and $\mathfrak{P}^{(2)}$, we have that
\begin{enumerate}
\item[(i)]
$\mathfrak{P}^{(2)}$ is a second-order path of length strictly grater than one containing a second-order echelon on a step different from zero.
\item[(ii)]
$\mathfrak{Q}^{(2)}$ is a second-order path.
\end{enumerate}

From (i) and taking into account the definition of $f^{(2)\flat}$ introduced in Proposition~\ref{PDPthExt}, we have that the value of the second-order path extension mapping $f^{(2)\flat}$ at $\mathfrak{P}^{(2)}$ is given by
$$
f_{s}^{(2)\flat}\left(
\mathfrak{P}^{(2)}
\right)
=
f_{s}^{(2)\flat}\left(
\mathfrak{P}^{(2)i,\bb{\mathfrak{P}^{(2)}}-1}
\right)
\circ^{1\mathbf{Pth}_{\boldsymbol{\mathcal{B}}^{(2)}}}_{\varphi(s)}
f_{s}^{(2)\flat}\left(
\mathfrak{P}^{(2)0,i-1}
\right).
$$

Since $(\mathfrak{Q}^{(2)} \circ_{s}^{1\mathbf{Pth}_{\boldsymbol{\mathcal{A}}^{(2)}}} \mathfrak{P}^{(2)i,\bb{\mathfrak{P}^{(2)}}-1}, s)$$\prec_{\mathbf{Pth}_{\boldsymbol{\mathcal{A}}^{(2)}}}$$(\mathfrak{Q}^{(2)} \circ_{s}^{1\mathbf{Pth}_{\boldsymbol{\mathcal{A}}^{(2)}}} \mathfrak{P}^{(2)}, s)$ we have, by induction, that
\begin{align*}
&\left\llbracket
f_{s}^{(2)\flat}\left(
\mathfrak{Q}^{(2)}
\circ_{s}^{1\mathbf{Pth}_{\boldsymbol{\mathcal{A}}^{(2)}}}
\mathfrak{P}^{(2)i,\bb{\mathfrak{P}^{(2)}}-1}
\right)
\right\rrbracket_{\varphi(s)}
\\
&\hspace{3cm}
=
\left\llbracket
f_{s}^{(2)\flat}\left(
\mathfrak{Q}^{(2)}
\right)
\right\rrbracket_{\varphi(s)}
\circ_{s}^{1\llbracket\mathbf{Pth}_{\boldsymbol{\mathcal{B}}^{(2)}}^{\mathbf{f}^{(2)}}\rrbracket}
\left\llbracket
f_{s}^{(2)\flat}\left(
\mathfrak{P}^{(2)i,\bb{\mathfrak{P}^{(2)}}-1}
\right)
\right\rrbracket_{\varphi(s)}.
\end{align*}

So, considering the foregoing, we can affirm that
\begin{flushleft}
$
\left\llbracket
f_{s}^{(2)\flat}\left(
\mathfrak{Q}^{(2)} \circ_{s}^{1\mathbf{Pth}_{\boldsymbol{\mathcal{A}}^{(2)}}} \mathfrak{P}^{(2)}
\right)
\right\rrbracket_{\varphi(s)}
$
\allowdisplaybreaks
\begin{align*}
&=
\left\llbracket
f_{s}^{(2)\flat}\left(
\mathfrak{Q}^{(2)} \circ_{s}^{1\mathbf{Pth}_{\boldsymbol{\mathcal{A}}^{(2)}}} \mathfrak{P}^{(2)i,\bb{\mathfrak{P}^{(2)}}-1}
\right)
\circ_{\varphi(s)}^{1\mathbf{Pth}_{\boldsymbol{\mathcal{B}}^{(2)}}}
f_{s}^{(2)\flat}\left(
\mathfrak{P}^{(2)0,i-1}
\right)
\right\rrbracket_{\varphi(s)}
\tag{1}
\\
&=
\left\llbracket
f_{s}^{(2)\flat}\left(
\mathfrak{Q}^{(2)} \circ_{s}^{1\mathbf{Pth}_{\boldsymbol{\mathcal{A}}^{(2)}}} \mathfrak{P}^{(2)i,\bb{\mathfrak{P}^{(2)}}-1}
\right)
\circ_{s}^{1\mathbf{Pth}_{\boldsymbol{\mathcal{B}}^{(2)}}^{\mathbf{f}^{(2)}}}
f_{s}^{(2)\flat}\left(
\mathfrak{P}^{(2)0,i-1}
\right)
\right\rrbracket_{\varphi(s)}
\tag{2}
\\
&=
\left\llbracket
f_{s}^{(2)\flat}\left(
\mathfrak{Q}^{(2)} \circ_{s}^{1\mathbf{Pth}_{\boldsymbol{\mathcal{A}}^{(2)}}} \mathfrak{P}^{(2)i,\bb{\mathfrak{P}^{(2)}}-1}
\right)
\right\rrbracket_{\varphi(s)}
\circ_{s}^{1\llbracket\mathbf{Pth}_{\boldsymbol{\mathcal{B}}^{(2)}}^{\mathbf{f}^{(2)}}\rrbracket}
\\
&\hspace{8cm}
\left\llbracket
f_{s}^{(2)\flat}\left(
\mathfrak{P}^{(2)0,i-1}
\right)
\right\rrbracket_{\varphi(s)}
\tag{3}
\\
&=
\left(
\left\llbracket
f_{s}^{(2)\flat}\left(
\mathfrak{Q}^{(2)}
\right)
\right\rrbracket_{\varphi(s)}
\circ_{s}^{1\llbracket\mathbf{Pth}_{\boldsymbol{\mathcal{B}}^{(2)}}^{\mathbf{f}^{(2)}}\rrbracket} 
\left\llbracket
f_{s}^{(2)\flat}\left(
\mathfrak{P}^{(2)i,\bb{\mathfrak{P}}-1}
\right)
\right\rrbracket_{\varphi(s)}
\right)
\\
&\hspace{6.5cm}
\circ_{s}^{1\llbracket\mathbf{Pth}_{\boldsymbol{\mathcal{B}}^{(2)}}^{\mathbf{f}^{(2)}}\rrbracket}
\left\llbracket
f_{s}^{(2)\flat}\left(
\mathfrak{P}^{(2)0,i-1}
\right)
\right\rrbracket_{\varphi(s)}
\tag{4}
\\
&=
\left\llbracket
f_{s}^{(2)\flat}\left(
\mathfrak{Q}^{(2)}
\right)
\right\rrbracket_{\varphi(s)}
\circ_{s}^{1\llbracket\mathbf{Pth}_{\boldsymbol{\mathcal{B}}^{(2)}}^{\mathbf{f}^{(2)}}\rrbracket} 
\\
&\hspace{2cm}
\left(
\left\llbracket
f_{s}^{(2)\flat}\left(
\mathfrak{P}^{(2)i,\bb{\mathfrak{P}}-1}
\right)
\right\rrbracket_{\varphi(s)}
\circ_{s}^{1\llbracket\mathbf{Pth}_{\boldsymbol{\mathcal{B}}^{(2)}}^{\mathbf{f}^{(2)}}\rrbracket}
\left\llbracket
f_{s}^{(2)\flat}\left(
\mathfrak{P}^{(2)0,i-1}
\right)
\right\rrbracket_{\varphi(s)}
\right)
\tag{5}
\\
&=
\left\llbracket
f_{s}^{(2)\flat}\left(
\mathfrak{Q}^{(2)}
\right)
\right\rrbracket_{\varphi(s)}
\circ_{s}^{1\llbracket\mathbf{Pth}_{\boldsymbol{\mathcal{B}}^{(2)}}^{\mathbf{f}^{(2)}}\rrbracket} 
\\
&\hspace{3.5cm}
\left\llbracket
f_{s}^{(2)\flat}\left(
\mathfrak{P}^{(2)i,\bb{\mathfrak{P}^{(2)}}-1}
\right)
\circ_{s}^{1\mathbf{Pth}_{\boldsymbol{\mathcal{B}}^{(2)}}^{\mathbf{f}^{(2)}}}
f_{s}^{(2)\flat}\left(
\mathfrak{P}^{(2)0,i-1}
\right)
\right\rrbracket_{\varphi(s)}
\tag{6}
\\
&=
\left\llbracket
f_{s}^{(2)\flat}\left(
\mathfrak{Q}^{(2)}
\right)
\right\rrbracket_{\varphi(s)}
\circ_{s}^{1\llbracket\mathbf{Pth}_{\boldsymbol{\mathcal{B}}^{(2)}}^{\mathbf{f}^{(2)}}\rrbracket}
\\
&\hspace{3.5cm}
\left\llbracket
f_{s}^{(2)\flat}\left(
\mathfrak{P}^{(2)i,\bb{\mathfrak{P}}-1}
\right)
\circ_{\varphi(s)}^{1\mathbf{Pth}_{\boldsymbol{\mathcal{B}}^{(1)}}}
f_{s}^{(2)\flat}\left(
\mathfrak{P}^{(2)0,i-1}
\right)
\right\rrbracket_{\varphi(s)}
\tag{7}
\\
&=
\left\llbracket
f_{s}^{(2)\flat}\left(
\mathfrak{Q}^{(2)}
\right)
\right\rrbracket_{\varphi(s)}
\circ_{s}^{1\llbracket\mathbf{Pth}_{\boldsymbol{\mathcal{B}}^{(2)}}^{\mathbf{f}^{(2)}}\rrbracket}
\\
&\hspace{3.5cm}
\left\llbracket
f_{s}^{(2)\flat}\left(
\mathfrak{P}^{(2)i,\bb{\mathfrak{P}^{(2)}}-1}
\circ_{s}^{1\mathbf{Pth}_{\boldsymbol{\mathcal{A}}^{(2)}}}
\mathfrak{P}^{(2)0,i-1}
\right)
\right\rrbracket_{\varphi(s)}
\tag{8}
\\
&=
\left\llbracket
f_{s}^{(2)\flat}\left(
\mathfrak{Q}^{(2)}
\right)
\right\rrbracket_{\varphi(s)}
\circ_{s}^{1\llbracket\mathbf{Pth}_{\boldsymbol{\mathcal{B}}^{(2)}}^{\mathbf{f}^{(2)}}\rrbracket} 
\left\llbracket
f_{s}^{(2)\flat}\left(
\mathfrak{P}^{(2)}
\right)
\right\rrbracket_{\varphi(s)}.
\tag{9}
\end{align*}
\end{flushleft}

The first equality unravels the definition of the second-order path extension mapping $f^{(2)\flat}$ introduced in Proposition~\ref{PDPthExt};
the second equality recovers the definition of the operation symbol $\circ_{s}^{1}$ in the partial $\Sigma^{\boldsymbol{\mathcal{A}}^{(2)}}$-algebra $\mathbf{Pth}_{\boldsymbol{\mathcal{B}}^{(2)}}^{\mathbf{f}^{(2)}}$ introduced in Proposition~\ref{PDPthBCatAlg};
the third equality recovers the definition of the operation symbol $\circ_{s}^{1}$ in the partial $\Sigma^{\boldsymbol{\mathcal{A}}^{(1)}}$-algebra $\llbracket\mathbf{Pth}_{\boldsymbol{\mathcal{B}}^{(2)}}^{\mathbf{f}^{(2)}}\rrbracket$ introduced in Proposition~\ref{PDQPthBCatAlg};
the fourth equality follows by Artinian induction;
the fifth equality follows from the fact that, according to Proposition~\ref{PDPthComp}, the interpretation of the partial operation $\circ_{s}^{1}$ is associative;
the sixth equality unravels the definition of the operation symbol $\circ_{s}^{1}$ in the partial $\Sigma^{\boldsymbol{\mathcal{A}}^{(2)}}$-algebra $\llbracket\mathbf{Pth}_{\boldsymbol{\mathcal{B}}^{(2)}}^{\mathbf{f}^{(2)}}\rrbracket$ introduced in Proposition~\ref{PDQPthBCatAlg};
The seventh equality unravels the definition of the operation symbol $\circ_{s}^{1}$ in the partial $\Sigma^{\boldsymbol{\mathcal{A}}^{(2)}}$-algebra $\mathbf{Pth}_{\boldsymbol{\mathcal{B}}^{(2)}}^{\mathbf{f}^{(2)}}$ introduced in Proposition~\ref{PDPthBCatAlg};
the eight equality recovers the definition of $f^{(2)\flat}$ since $\mathfrak{P}^{(2)}$ is a second-order path of length strictly greater than one containing a second-order echelon on its first step;
finally, the last equality recovers the definition of $\mathfrak{P}^{(2)}$.

If~(1.2.2), i.e, $i \neq 0$ and $i \in [\bb{\mathfrak{P}^{(2)}}, \bb{\mathfrak{Q}^{(2)} \circ_{s}^{1\mathbf{Pth}_{\boldsymbol{\mathcal{A}}^{(2)}} \mathfrak{P}^{(2)}}-1}]$, then $\mathfrak{Q}^{(2)}$ is not an $(2,[1])$-identity path containing a second-order echelon, whilst $\mathfrak{P}^{(2)}$ is a second-order echelonless path.

We will distinguish three cases according to whether (1.2.2.1) $\mathfrak{Q}^{(2)}$ is a second-order echelon; (1.2.2.2) $\mathfrak{Q}^{(2)}$ is a second-order path of length strictly greater than one containing a second-order echelon of its first step or (1.2.2.3) $\mathfrak{Q}^{(2)}$ is a second-order path of length strictly greater than one containing a second-order echelon on a step different from zero. These cases can be proved using a similar argument to those three cases presented above. We leave the details for the interested reader.

This finishes the case $i > 0$.

This completes case $(1)$.

If~(2), i.e., if $\mathfrak{Q}^{(2)} \circ_{s}^{1\mathbf{Pth}_{\boldsymbol{\mathcal{A}}^{(2)}}} \mathfrak{P}^{(2)}$ is an echelonless second-order path, it could be the case that $(2.1)$ $\mathfrak{Q}^{(2)} \circ_{s}^{1\mathbf{Pth}_{\boldsymbol{\mathcal{A}}^{(2)}}} \mathfrak{P}^{(2)}$ is an echelonless second-order path that is not head-constant, or $(2.2)$ $\mathfrak{Q}^{(2)} \circ_{s}^{1\mathbf{Pth}_{\boldsymbol{\mathcal{A}}^{(2)}}} \mathfrak{P}^{(2)}$ is a head-constant echelonless second-order path that is not coherent, or $(2.3)$ $\mathfrak{Q}^{(2)} \circ_{s}^{1\mathbf{Pth}_{\boldsymbol{\mathcal{A}}^{(2)}}} \mathfrak{P}^{(2)}$ is a coherent head-constant echelonless second-order path.

If $(2.1)$, let $i \in \bb{\mathfrak{Q}^{(2)} \circ_{s}^{1\mathbf{Pth}_{\boldsymbol{\mathcal{A}}^{(2)}}} \mathfrak{P}^{(2)}}$ be the greatest index for which $(\mathfrak{Q}^{(2)} \circ_{s}^{1\mathbf{Pth}_{\boldsymbol{\mathcal{A}}^{(2)}}} \mathfrak{P}^{(2)})^{0,i}$ is a head- constant second-order path. Since $\bb{\mathfrak{Q}^{(2)} \circ_{s}^{1\mathbf{Pth}_{\boldsymbol{\mathcal{A}}^{(2)}}} \mathfrak{P}^{(2)}} = \bb{\mathfrak{Q}^{(2)}} + \bb{\mathfrak{P}^{(2)}}$, we have that either $(2.1.1)$ $i \in \bb{\mathfrak{P}^{(2)}}-1$, $(2.1.2)$ $i = \bb{\mathfrak{P}^{(2)}}-1$, or $(2.1.3)$ $i \in [\bb{\mathfrak{P}^{(2)}}, \bb{\mathfrak{Q}^{(2)} \circ_{s}^{1\mathbf{Pth}_{\boldsymbol{\mathcal{A}}^{(2)}}} \mathfrak{P}^{(2)}}-1]$.

If $(2.1.1)$, i.e., if we find ourselves in the case where $\mathfrak{Q}^{(2)} \circ_{s}^{1\mathbf{Pth}_{\boldsymbol{\mathcal{A}}^{(2)}}} \mathfrak{P}^{(2)}$ is an echelonless second-order path that is not head-constant and $i \in \bb{\mathfrak{P}^{(2)}}-1$ is the greatest index for which $(\mathfrak{Q}^{(2)} \circ_{s}^{1\mathbf{Pth}_{\boldsymbol{\mathcal{A}}^{(2)}}} \mathfrak{P}^{(2)})^{0,i}$ is head-constant then, regarding the second-order paths $\mathfrak{Q}^{(2)}$ and $\mathfrak{P}^{(2)}$, we have that
\begin{enumerate}
\item[(i)]
$\mathfrak{P}^{(2)}$ is an echelonless second-order path that is not head-constant and $i \in \bb{\mathfrak{P}^{(2)}}-1$ is the greatest index for which $\mathfrak{P}^{(2)0,i}$ is head-constant.
\item[(ii)]
$\mathfrak{Q}^{(2)}$ is an echelonless second-order path.
\end{enumerate}

From (i) and taking into account Proposition~\ref{PDPthExt}, we have that the value of the path extension mapping of $\mathbf{f}^{(2)}$ at $\mathfrak{P}^{(2)}$ is given by
$$
f_{s}^{(2)\flat}\left(
\mathfrak{P}^{(2)}
\right)
=
f_{s}^{(2)\flat}\left(
\mathfrak{P}^{(2)i+1,\bb{\mathfrak{P}^{(2)}}-1}
\right)
\circ_{\varphi(s)}^{1\mathbf{Pth}_{\boldsymbol{\mathcal{B}}^{(2)}}}
f_{s}^{(2)\flat}\left(
\mathfrak{P}^{(2)0,i}
\right).
$$

Since $(\mathfrak{Q}^{(2)} \circ_{s}^{1\mathbf{Pth}_{\boldsymbol{\mathcal{A}}^{(2)}}} \mathfrak{P}^{(2)i+1,\bb{\mathfrak{P}^{(2)}}-1}, s)$$\prec_{\mathbf{Pth}_{\boldsymbol{\mathcal{A}}^{(2)}}}$$(\mathfrak{Q}^{(2)} \circ_{s}^{1\mathbf{Pth}_{\boldsymbol{\mathcal{A}}^{(2)}}} \mathfrak{P}^{(2)}, s)$ we have, by induction, that
\begin{align*}
&\left\llbracket
f_{s}^{(2)\flat}\left(
\mathfrak{Q}^{(2)}
\circ_{s}^{1\mathbf{Pth}_{\boldsymbol{\mathcal{A}}^{(2)}}}
\mathfrak{P}^{(2)i+1,\bb{\mathfrak{P}^{(2)}}-1}
\right)
\right\rrbracket_{\varphi(s)}
\\
&\hspace{3cm}
=
\left\llbracket
f_{s}^{(2)\flat}\left(
\mathfrak{Q}^{(2)}
\right)
\right\rrbracket_{\varphi(s)}
\circ_{s}^{1\llbracket\mathbf{Pth}_{\boldsymbol{\mathcal{B}}^{(2)}}^{\mathbf{f}^{(2)}}\rrbracket}
\left\llbracket
f_{s}^{(2)\flat}\left(
\mathfrak{P}^{(2)i+1,\bb{\mathfrak{P}^{(2)}}-1}
\right)
\right\rrbracket_{\varphi(s)}.
\end{align*}

So, considering the foregoing, we can affirm that
\begin{flushleft}
$
\left\llbracket
f_{s}^{(2)\flat}\left(
\mathfrak{Q}^{(2)} \circ_{s}^{1\mathbf{Pth}_{\boldsymbol{\mathcal{A}}^{(2)}}} \mathfrak{P}^{(2)}
\right)
\right\rrbracket_{\varphi(s)}
$
\allowdisplaybreaks
\begin{align*}
&=
\left\llbracket
f_{s}^{(2)\flat}\left(
\mathfrak{Q}^{(2)} \circ_{s}^{1\mathbf{Pth}_{\boldsymbol{\mathcal{A}}^{(2)}}} \mathfrak{P}^{(2)i+1,\bb{\mathfrak{P}^{(2)}}-1}
\right)
\circ_{\varphi(s)}^{1\mathbf{Pth}_{\boldsymbol{\mathcal{B}}^{(2)}}}
f_{s}^{(2)\flat}\left(
\mathfrak{P}^{(2)0,i}
\right)
\right\rrbracket_{\varphi(s)}
\tag{1}
\\
&=
\left\llbracket
f_{s}^{(2)\flat}\left(
\mathfrak{Q}^{(2)} \circ_{s}^{1\mathbf{Pth}_{\boldsymbol{\mathcal{A}}^{(2)}}} \mathfrak{P}^{(2)i+1,\bb{\mathfrak{P}^{(2)}}-1}
\right)
\circ_{s}^{1\mathbf{Pth}_{\boldsymbol{\mathcal{B}}^{(2)}}^{\mathbf{f}^{(2)}}}
f_{s}^{(2)\flat}\left(
\mathfrak{P}^{(2)0,i}
\right)
\right\rrbracket_{\varphi(s)}
\tag{2}
\\
&=
\left\llbracket
f_{s}^{(2)\flat}\left(
\mathfrak{Q}^{(2)} \circ_{s}^{1\mathbf{Pth}_{\boldsymbol{\mathcal{A}}^{(2)}}} \mathfrak{P}^{(2)i+1,\bb{\mathfrak{P}^{(2)}}-1}
\right)
\right\rrbracket_{\varphi(s)}
\circ_{s}^{1\llbracket\mathbf{Pth}_{\boldsymbol{\mathcal{B}}^{(2)}}^{\mathbf{f}^{(2)}}\rrbracket}
\\
&\hspace{8cm}
\left\llbracket
f_{s}^{(2)\flat}\left(
\mathfrak{P}^{(2)0,i}
\right)
\right\rrbracket_{\varphi(s)}
\tag{3}
\\
&=
\left(
\left\llbracket
f_{s}^{(2)\flat}\left(
\mathfrak{Q}^{(2)}
\right)
\right\rrbracket_{\varphi(s)}
\circ_{s}^{1\llbracket\mathbf{Pth}_{\boldsymbol{\mathcal{B}}^{(2)}}^{\mathbf{f}^{(2)}}\rrbracket} 
\left\llbracket
f_{s}^{(2)\flat}\left(
\mathfrak{P}^{(2)i+1,\bb{\mathfrak{P}}-1}
\right)
\right\rrbracket_{\varphi(s)}
\right)
\\
&\hspace{6.5cm}
\circ_{s}^{1\llbracket\mathbf{Pth}_{\boldsymbol{\mathcal{B}}^{(2)}}^{\mathbf{f}^{(2)}}\rrbracket}
\left\llbracket
f_{s}^{(2)\flat}\left(
\mathfrak{P}^{(2)0,i}
\right)
\right\rrbracket_{\varphi(s)}
\tag{4}
\\
&=
\left\llbracket
f_{s}^{(2)\flat}\left(
\mathfrak{Q}^{(2)}
\right)
\right\rrbracket_{\varphi(s)}
\circ_{s}^{1\llbracket\mathbf{Pth}_{\boldsymbol{\mathcal{B}}^{(2)}}^{\mathbf{f}^{(2)}}\rrbracket} 
\\
&\hspace{2cm}
\left(
\left\llbracket
f_{s}^{(2)\flat}\left(
\mathfrak{P}^{(2)i+1,\bb{\mathfrak{P}}-1}
\right)
\right\rrbracket_{\varphi(s)}
\circ_{s}^{1\llbracket\mathbf{Pth}_{\boldsymbol{\mathcal{B}}^{(2)}}^{\mathbf{f}^{(2)}}\rrbracket}
\left\llbracket
f_{s}^{(2)\flat}\left(
\mathfrak{P}^{(2)0,i}
\right)
\right\rrbracket_{\varphi(s)}
\right)
\tag{5}
\\
&=
\left\llbracket
f_{s}^{(2)\flat}\left(
\mathfrak{Q}^{(2)}
\right)
\right\rrbracket_{\varphi(s)}
\circ_{s}^{1\llbracket\mathbf{Pth}_{\boldsymbol{\mathcal{B}}^{(2)}}^{\mathbf{f}^{(2)}}\rrbracket} 
\\
&\hspace{3.5cm}
\left\llbracket
f_{s}^{(2)\flat}\left(
\mathfrak{P}^{(2)i+1,\bb{\mathfrak{P}^{(2)}}-1}
\right)
\circ_{s}^{1\mathbf{Pth}_{\boldsymbol{\mathcal{B}}^{(2)}}^{\mathbf{f}^{(2)}}}
f_{s}^{(2)\flat}\left(
\mathfrak{P}^{(2)0,i}
\right)
\right\rrbracket_{\varphi(s)}
\tag{6}
\\
&=
\left\llbracket
f_{s}^{(2)\flat}\left(
\mathfrak{Q}^{(2)}
\right)
\right\rrbracket_{\varphi(s)}
\circ_{s}^{1\llbracket\mathbf{Pth}_{\boldsymbol{\mathcal{B}}^{(2)}}^{\mathbf{f}^{(2)}}\rrbracket}
\\
&\hspace{3.5cm}
\left\llbracket
f_{s}^{(2)\flat}\left(
\mathfrak{P}^{(2)i+1,\bb{\mathfrak{P}}-1}
\right)
\circ_{\varphi(s)}^{1\mathbf{Pth}_{\boldsymbol{\mathcal{B}}^{(2)}}}
f_{s}^{(2)\flat}\left(
\mathfrak{P}^{(2)0,i}
\right)
\right\rrbracket_{\varphi(s)}
\tag{7}
\\
&=
\left\llbracket
f_{s}^{(2)\flat}\left(
\mathfrak{Q}^{(2)}
\right)
\right\rrbracket_{\varphi(s)}
\circ_{s}^{1\llbracket\mathbf{Pth}_{\boldsymbol{\mathcal{B}}^{(2)}}^{\mathbf{f}^{(2)}}\rrbracket}
\\
&\hspace{3.5cm}
\left\llbracket
f_{s}^{(2)\flat}\left(
\mathfrak{P}^{(2)i+1,\bb{\mathfrak{P}^{(2)}}-1}
\circ_{s}^{1\mathbf{Pth}_{\boldsymbol{\mathcal{A}}^{(2)}}}
\mathfrak{P}^{(2)0,i}
\right)
\right\rrbracket_{\varphi(s)}
\tag{8}
\\
&=
\left\llbracket
f_{s}^{(2)\flat}\left(
\mathfrak{Q}^{(2)}
\right)
\right\rrbracket_{\varphi(s)}
\circ_{s}^{1\llbracket\mathbf{Pth}_{\boldsymbol{\mathcal{B}}^{(2)}}^{\mathbf{f}^{(2)}}\rrbracket} 
\left\llbracket
f_{s}^{(2)\flat}\left(
\mathfrak{P}^{(2)}
\right)
\right\rrbracket_{\varphi(s)}.
\tag{9}
\end{align*}
\end{flushleft}

The first equality unravels the definition of the second-order path extension mapping $f^{(2)\flat}$ introduced in Proposition~\ref{PDPthExt};
the second equality recovers the definition of the operation symbol $\circ_{s}^{1}$ in the partial $\Sigma^{\boldsymbol{\mathcal{A}}^{(2)}}$-algebra $\mathbf{Pth}_{\boldsymbol{\mathcal{B}}^{(2)}}^{\mathbf{f}^{(2)}}$ introduced in Proposition~\ref{PDPthBCatAlg};
the third equality recovers the definition of the operation symbol $\circ_{s}^{1}$ in the partial $\Sigma^{\boldsymbol{\mathcal{A}}^{(1)}}$-algebra $\llbracket\mathbf{Pth}_{\boldsymbol{\mathcal{B}}^{(2)}}^{\mathbf{f}^{(2)}}\rrbracket$ introduced in Proposition~\ref{PDQPthBCatAlg};
the fourth equality follows by Artinian induction;
the fifth equality follows from the fact that, according to Proposition~\ref{PDPthComp}, the interpretation of the partial operation $\circ_{s}^{1}$ is associative;
the sixth equality unravels the definition of the operation symbol $\circ_{s}^{1}$ in the partial $\Sigma^{\boldsymbol{\mathcal{A}}^{(2)}}$-algebra $\llbracket\mathbf{Pth}_{\boldsymbol{\mathcal{B}}^{(2)}}^{\mathbf{f}^{(2)}}\rrbracket$ introduced in Proposition~\ref{PDQPthBCatAlg};
the seventh equality unravels the definition of the operation symbol $\circ_{s}^{1}$ in the partial $\Sigma^{\boldsymbol{\mathcal{A}}^{(2)}}$-algebra $\mathbf{Pth}_{\boldsymbol{\mathcal{B}}^{(2)}}^{\mathbf{f}^{(2)}}$ introduced in Proposition~\ref{PDPthBCatAlg};
the eight equality recovers the definition of the path extension mapping $f^{(2)\flat}$ introduced in Proposition~\ref{PPthExt};
finally, the last equality recovers the definition of $\mathfrak{P}^{(2)}$.

The case $i \in \bb{\mathfrak{P}^{(2)}}-1$ follows.

If $(2.1.2)$, i.e., if we find ourselves in the case where $\mathfrak{Q}^{(2)} \circ_{s}^{1\mathbf{Pth}_{\boldsymbol{\mathcal{A}}^{(2)}}} \mathfrak{P}^{(2)}$ is an echelonless second-order path that is not head-constant and $i = \bb{\mathfrak{P}^{(2)}}-1$ is the greatest index for which $(\mathfrak{Q}^{(2)} \circ_{s}^{1\mathbf{Pth}_{\boldsymbol{\mathcal{A}}^{(2)}}} \mathfrak{P}^{(2)})^{0,i}$ is head-constant then, regarding the second-order paths $\mathfrak{Q}^{(2)}$ and $\mathfrak{P}^{(2)}$, we have that
\begin{enumerate}
\item[(i)]
$\mathfrak{P}^{(2)}$ is a head-constant echelonless second-order path.
\item[(ii)]
$\mathfrak{Q}^{(2)}$ is an echelonless second-order path.
\end{enumerate}

So, considering the foregoing, we can affirm that
\begin{flushleft}
$
\left\llbracket
f_{s}^{(2)\flat}\left(
\mathfrak{Q}^{(2)} \circ_{s}^{1\mathbf{Pth}_{\boldsymbol{\mathcal{A}}^{(2)}}} \mathfrak{P}^{(2)}
\right)
\right\rrbracket_{\varphi(s)}
$
\allowdisplaybreaks
\begin{align*}
&=
\left\llbracket
f_{s}^{(2)\flat}\left(
\mathfrak{Q}^{(2)}
\right)
\circ_{\varphi(s)}^{1\mathbf{Pth}_{\boldsymbol{\mathcal{B}}^{(2)}}}
f_{s}^{(2)\flat}\left(
\mathfrak{P}^{(2)}
\right)
\right\rrbracket_{\varphi(s)}
\tag{1}
\\
&=
\left\llbracket
f_{s}^{(2)\flat}\left(
\mathfrak{Q}^{(2)}
\right)
\circ_{s}^{1\mathbf{Pth}_{\boldsymbol{\mathcal{B}}^{(2)}}^{\mathbf{f}^{(2)}}}
f_{s}^{(2)\flat}\left(
\mathfrak{P}^{(2)}
\right)
\right\rrbracket_{\varphi(s)}
\tag{2}
\\
&=
\left\llbracket
f_{s}^{(2)\flat}\left(
\mathfrak{Q}^{(2)}
\right)
\right\rrbracket_{\varphi(s)}
\circ_{s}^{1\llbracket\mathbf{Pth}_{\boldsymbol{\mathcal{B}}^{(2)}}^{\mathbf{f}^{(2)}}\rrbracket}
\left\llbracket
f_{s}^{(2)\flat}\left(
\mathfrak{P}^{(2)}
\right)
\right\rrbracket_{\varphi(s)}
\tag{3}
\end{align*}
\end{flushleft}

The first equality unravels the definition of the path extension mapping $f^{(2)\flat}$ introduced in Proposition~\ref{PDPthExt};
the second equality recovers the definition of the operation symbol $\circ_{s}^{1}$ in the partial $\Sigma^{\boldsymbol{\mathcal{A}}^{(2)}}$-algebra $\mathbf{Pth}_{\boldsymbol{\mathcal{B}}^{(2)}}^{\mathbf{f}^{(2)}}$ introduced in Proposition~\ref{PDPthBCatAlg};
finally, the last equality recovers the definition of the operation symbol $\circ_{s}^{1}$ in the partial $\Sigma^{\boldsymbol{\mathcal{A}}^{(1)}}$-algebra $\llbracket\mathbf{Pth}_{\boldsymbol{\mathcal{B}}^{(2)}}^{\mathbf{f}^{(2)}}\rrbracket$ introduced in Proposition~\ref{PDQPthBCatAlg}.

The case $i = \bb{\mathfrak{P}^{(2)}}-1$ follows.

If $(2.1.3)$, i.e., if we find ourselves in the case where $\mathfrak{Q}^{(2)} \circ_{s}^{1\mathbf{Pth}_{\boldsymbol{\mathcal{A}}^{(2)}}} \mathfrak{P}^{(2)}$ is an echelonless second-order path that is not head-constant and $i \in [\bb{\mathfrak{P}^{(2)}}, \bb{\mathfrak{Q}^{(2)} \circ_{s}^{1\mathbf{Pth}_{\boldsymbol{\mathcal{A}}^{(2)}}} \mathfrak{P}^{(2)}}-1]$ is the greatest index for which $(\mathfrak{Q}^{(2)} \circ_{s}^{1\mathbf{Pth}_{\boldsymbol{\mathcal{A}}^{(2)}}} \mathfrak{P}^{(2)})^{0,i}$ is head-constant then, regarding the second-order paths $\mathfrak{Q}^{(2)}$ and $\mathfrak{P}^{(2)}$, we have that
\begin{enumerate}
\item[(i)]
$\mathfrak{P}^{(2)}$ is a head-constant echelonless second-order path.
\item[(ii)]
$\mathfrak{Q}^{(2)}$ is an echelonless second-order path that is not head-constant and $j = i - \bb{\mathfrak{P}^{(2)}} \in \bb{\mathfrak{Q}^{(2)}}-1$ is the greatest index for which $\mathfrak{Q}^{(2)0,j}$ is head-constant.
\end{enumerate}

From (ii) and taking into account Proposition~\ref{PDPthExt}, we have that the value of the path extension mapping of $\mathbf{f}^{(2)}$ at $\mathfrak{Q}^{(2)}$ is given by
$$
f_{s}^{(2)\flat}\left(
\mathfrak{Q}^{(2)}
\right)
=
f_{s}^{(2)\flat}\left(
\mathfrak{Q}^{(2)j+1,\bb{\mathfrak{Q}^{(2)}}-1}
\right)
\circ_{\varphi(s)}^{1\mathbf{Pth}_{\boldsymbol{\mathcal{B}}^{(2)}}}
f_{s}^{(2)\flat}\left(
\mathfrak{Q}^{(2)0,j}
\right).
$$

Since $(\mathfrak{Q}^{(2)0, j} \circ_{s}^{1\mathbf{Pth}_{\boldsymbol{\mathcal{A}}^{(2)}}} \mathfrak{P}^{(2)}, s)$$\prec_{\mathbf{Pth}_{\boldsymbol{\mathcal{A}}^{(2)}}}$$(\mathfrak{Q}^{(2)} \circ_{s}^{1\mathbf{Pth}_{\boldsymbol{\mathcal{A}}^{(2)}}} \mathfrak{P}^{(2)}, s)$ we have, by induction, that
\begin{align*}
&\left\llbracket
f_{s}^{(2)\flat}\left(
\mathfrak{Q}^{(2)0, j}
\circ_{s}^{1\mathbf{Pth}_{\boldsymbol{\mathcal{A}}^{(2)}}}
\mathfrak{P}^{(2)}
\right)
\right\rrbracket_{\varphi(s)}
\\
&\hspace{3cm}
=
\left\llbracket
f_{s}^{(2)\flat}\left(
\mathfrak{Q}^{(2)0, j}
\right)
\right\rrbracket_{\varphi(s)}
\circ_{s}^{1\llbracket\mathbf{Pth}_{\boldsymbol{\mathcal{B}}^{(2)}}^{\mathbf{f}^{(2)}}\rrbracket}
\left\llbracket
f_{s}^{(2)\flat}\left(
\mathfrak{P}^{(2)}
\right)
\right\rrbracket_{\varphi(s)}.
\end{align*}

So, considering the foregoing, we can affirm that
\begin{flushleft}
$
\left\llbracket
f_{s}^{(2)\flat}\left(
\mathfrak{Q}^{(2)} \circ_{s}^{1\mathbf{Pth}_{\boldsymbol{\mathcal{A}}^{(2)}}} \mathfrak{P}^{(2)}
\right)
\right\rrbracket_{\varphi(s)}
$
\allowdisplaybreaks
\begin{align*}
&=
\left\llbracket
f_{s}^{(2)\flat}\left(
\mathfrak{Q}^{(2)j+1, \bb{\mathfrak{Q}^{(2)}}-1}
\right)
\circ_{\varphi(s)}^{1\mathbf{Pth}_{\boldsymbol{\mathcal{B}}^{(2)}}}
f_{s}^{(2)\flat}\left(
\mathfrak{Q}^{(2)0, j}
\circ_{s}^{1\mathbf{Pth}_{\boldsymbol{\mathcal{A}}^{(2)}}}
\mathfrak{P}^{(2)}
\right)
\right\rrbracket_{\varphi(s)}
\tag{1}
\\
&=
\left\llbracket
f_{s}^{(2)\flat}\left(
\mathfrak{Q}^{(2)j+1, \bb{\mathfrak{Q}^{(2)}}-1}
\right)
\circ_{s}^{1\mathbf{Pth}_{\boldsymbol{\mathcal{B}}^{(2)}}^{\mathbf{f}^{(2)}}}
f_{s}^{(2)\flat}\left(
\mathfrak{Q}^{(2)0, j}
\circ_{s}^{1\mathbf{Pth}_{\boldsymbol{\mathcal{A}}^{(2)}}}
\mathfrak{P}^{(2)}
\right)
\right\rrbracket_{\varphi(s)}
\tag{2}
\\
&=
\left\llbracket
f_{s}^{(2)\flat}\left(
\mathfrak{Q}^{(2)j+1, \bb{\mathfrak{Q}^{(2)}}-1}
\right)
\right\rrbracket_{\varphi(s)}
\circ_{s}^{1\llbracket\mathbf{Pth}_{\boldsymbol{\mathcal{B}}^{(2)}}^{\mathbf{f}^{(2)}}\rrbracket}
\\
&\hspace{6cm}
\left\llbracket
f_{s}^{(2)\flat}\left(
\mathfrak{Q}^{(2)0, j}
\circ_{s}^{1\mathbf{Pth}_{\boldsymbol{\mathcal{A}}^{(2)}}}
\mathfrak{P}^{(2)}
\right)
\right\rrbracket_{\varphi(s)}
\tag{3}
\\
&=
\left\llbracket
f_{s}^{(2)\flat}\left(
\mathfrak{Q}^{(2)j+1, \bb{\mathfrak{Q}^{(2)}}-1}
\right)
\right\rrbracket_{\varphi(s)}
\circ_{s}^{1\llbracket\mathbf{Pth}_{\boldsymbol{\mathcal{B}}^{(2)}}^{\mathbf{f}^{(2)}}\rrbracket}
\\
&\hspace{3cm}
\left(
\left\llbracket
f_{s}^{(2)\flat}\left(
\mathfrak{Q}^{(2)0, j}
\right)
\right\rrbracket_{\varphi(s)}
\circ_{s}^{1\llbracket\mathbf{Pth}_{\boldsymbol{\mathcal{B}}^{(2)}}^{\mathbf{f}^{(2)}}\rrbracket}
\left\llbracket
f_{s}^{(2)\flat}\left(
\mathfrak{P}^{(2)}
\right)
\right\rrbracket_{\varphi(s)}
\right)
\tag{4}
\\
&=
\left(
\left\llbracket
f_{s}^{(2)\flat}\left(
\mathfrak{Q}^{(2)j+1, \bb{\mathfrak{Q}^{(2)}}-1}
\right)
\right\rrbracket_{\varphi(s)}
\circ_{s}^{1\llbracket\mathbf{Pth}_{\boldsymbol{\mathcal{B}}^{(2)}}^{\mathbf{f}^{(2)}}\rrbracket}
\left\llbracket
f_{s}^{(2)\flat}\left(
\mathfrak{Q}^{(2)0, j}
\right)
\right\rrbracket_{\varphi(s)}
\right)
\\
&\hspace{7cm}
\circ_{s}^{1\llbracket\mathbf{Pth}_{\boldsymbol{\mathcal{B}}^{(2)}}^{\mathbf{f}^{(2)}}\rrbracket}
\left\llbracket
f_{s}^{(2)\flat}\left(
\mathfrak{P}^{(2)}
\right)
\right\rrbracket_{\varphi(s)}
\tag{5}
\\
&=
\left\llbracket
f_{s}^{(2)\flat}\left(
\mathfrak{Q}^{(2)j+1, \bb{\mathfrak{Q}^{(2)}}-1}
\right)
\circ_{s}^{1\mathbf{Pth}_{\boldsymbol{\mathcal{B}}^{(2)}}^{\mathbf{f}^{(2)}}}
f_{s}^{(2)\flat}\left(
\mathfrak{Q}^{(2)0, j}
\right)
\right\rrbracket_{\varphi(s)}
\\
&\hspace{7cm}
\circ_{s}^{1\llbracket\mathbf{Pth}_{\boldsymbol{\mathcal{B}}^{(2)}}^{\mathbf{f}^{(2)}}\rrbracket}
\left\llbracket
f_{s}^{(2)\flat}\left(
\mathfrak{P}^{(2)}
\right)
\right\rrbracket_{\varphi(s)}
\tag{6}
\\
&=
\left\llbracket
f_{s}^{(2)\flat}\left(
\mathfrak{Q}^{(2)j+1, \bb{\mathfrak{Q}^{(2)}}-1}
\right)
\circ_{\varphi(s)}^{1\mathbf{Pth}_{\boldsymbol{\mathcal{B}}^{(2)}}}
f_{s}^{(2)\flat}\left(
\mathfrak{Q}^{(2)0, j}
\right)
\right\rrbracket_{\varphi(s)}
\\
&\hspace{7cm}
\circ_{s}^{1\llbracket\mathbf{Pth}_{\boldsymbol{\mathcal{B}}^{(2)}}^{\mathbf{f}^{(2)}}\rrbracket}
\left\llbracket
f_{s}^{(2)\flat}\left(
\mathfrak{P}^{(2)}
\right)
\right\rrbracket_{\varphi(s)}
\tag{7}
\\
&=
\left\llbracket
f_{s}^{(2)\flat}\left(
\mathfrak{Q}^{(2)j+1, \bb{\mathfrak{Q}^{(2)}}-1}
\circ_{s}^{1\mathbf{Pth}_{\boldsymbol{\mathcal{A}}^{(2)}}}
\mathfrak{Q}^{(2)0, j}
\right)
\right\rrbracket_{\varphi(s)}
\\
&\hspace{7cm}
\circ_{s}^{1\llbracket\mathbf{Pth}_{\boldsymbol{\mathcal{B}}^{(2)}}^{\mathbf{f}^{(2)}}\rrbracket}
\left\llbracket
f_{s}^{(2)\flat}\left(
\mathfrak{P}^{(2)}
\right)
\right\rrbracket_{\varphi(s)}
\tag{8}
\\
&=
\left\llbracket
f_{s}^{(2)\flat}\left(
\mathfrak{Q}^{(2)}
\right)
\right\rrbracket_{\varphi(s)}
\circ_{s}^{1\llbracket\mathbf{Pth}_{\boldsymbol{\mathcal{B}}^{(2)}}^{\mathbf{f}^{(2)}}\rrbracket} 
\left\llbracket
f_{s}^{(2)\flat}\left(
\mathfrak{P}^{(2)}
\right)
\right\rrbracket_{\varphi(s)}.
\tag{9}
\end{align*}
\end{flushleft}

The first equality unravels the definition of the second-order path extension mapping $f^{(2)\flat}$ introduced in Proposition~\ref{PDPthExt};
the second equality recovers the definition of the operation symbol $\circ_{s}^{1}$ in the partial $\Sigma^{\boldsymbol{\mathcal{A}}^{(2)}}$-algebra $\mathbf{Pth}_{\boldsymbol{\mathcal{B}}^{(2)}}^{\mathbf{f}^{(2)}}$ introduced in Proposition~\ref{PDPthBCatAlg};
the third equality recovers the definition of the operation symbol $\circ_{s}^{1}$ in the partial $\Sigma^{\boldsymbol{\mathcal{A}}^{(1)}}$-algebra $\llbracket\mathbf{Pth}_{\boldsymbol{\mathcal{B}}^{(2)}}^{\mathbf{f}^{(2)}}\rrbracket$ introduced in Proposition~\ref{PDQPthBCatAlg};
the fourth equality follows by Artinian induction;
the fifth equality follows from the fact that, according to Proposition~\ref{PDPthComp}, the interpretation of the partial operation $\circ_{s}^{1}$ is associative;
the sixth equality unravels the definition of the operation symbol $\circ_{s}^{1}$ in the partial $\Sigma^{\boldsymbol{\mathcal{A}}^{(2)}}$-algebra $\llbracket\mathbf{Pth}_{\boldsymbol{\mathcal{B}}^{(2)}}^{\mathbf{f}^{(2)}}\rrbracket$ introduced in Proposition~\ref{PDQPthBCatAlg};
the seventh equality unravels the definition of the operation symbol $\circ_{s}^{1}$ in the partial $\Sigma^{\boldsymbol{\mathcal{A}}^{(2)}}$-algebra $\mathbf{Pth}_{\boldsymbol{\mathcal{B}}^{(2)}}^{\mathbf{f}^{(2)}}$ introduced in Proposition~\ref{PDPthBCatAlg};
the eight equality recovers the definition of the path extension mapping $f^{(2)\flat}$ introduced in Proposition~\ref{PPthExt};
finally, the last equality recovers the definition of $\mathfrak{Q}^{(2)}$.

The case $i \in [\bb{\mathfrak{P}^{(2)}}, \bb{\mathfrak{Q}^{(2)} \circ_{s}^{1\mathbf{Pth}_{\boldsymbol{\mathcal{A}}^{(2)}}} \mathfrak{P}^{(2)}}-1]$ follows.

This completes the case $(2.1)$.

If $(2.2)$, i.e., if $\mathfrak{Q}^{(2)} \circ_{s}^{1\mathbf{Pth}_{\boldsymbol{\mathcal{A}}^{(2)}}} \mathfrak{P}^{(2)}$ is a head-constant echelonless second-order path that is not coherent, then let $i \in \bb{\mathfrak{Q}^{(2)} \circ_{s}^{1\mathbf{Pth}_{\boldsymbol{\mathcal{A}}^{(2)}}} \mathfrak{P}^{(2)}}$ be the greatest index for which $(\bb{\mathfrak{Q}^{(2)} \circ_{s}^{1\mathbf{Pth}_{\boldsymbol{\mathcal{A}}^{(2)}}} \mathfrak{P}^{(2)}})^{0,i}$ is a coherent head-constant echelonless second-order path. Since $\bb{\mathfrak{Q}^{(2)} \circ_{s}^{1\mathbf{Pth}_{\boldsymbol{\mathcal{A}}^{(2)}}} \mathfrak{P}^{(2)}} = \bb{\mathfrak{Q}^{(2)}} + \bb{\mathfrak{P}^{(2)}}$, it follows that either $(2.1.1)$ $i \in \bb{\mathfrak{P}^{(2)}}-1$, $(2.1.2)$ $i = \bb{\mathfrak{P}^{(2)}}-1$, or $(2.1.3)$ $i \in [\bb{\mathfrak{P}^{(2)}}, \bb{\mathfrak{Q}^{(2)} \circ_{s}^{1\mathbf{Pth}_{\boldsymbol{\mathcal{A}}^{(2)}}} \mathfrak{P}^{(2)}}-1]$.

If $(2.2.1)$, i.e., if we find ourselves in the case where $\mathfrak{Q}^{(2)} \circ_{s}^{1\mathbf{Pth}_{\boldsymbol{\mathcal{A}}^{(2)}}} \mathfrak{P}^{(2)}$ is a head-constant echelonless second-order path that is not coherent and $i \in \bb{\mathfrak{P}^{(2)}}-1$ is the greatest index for which $(\mathfrak{Q}^{(2)} \circ_{s}^{1\mathbf{Pth}_{\boldsymbol{\mathcal{A}}^{(2)}}} \mathfrak{P}^{(2)})^{0,i}$ is coherent then, regarding the second-order paths $\mathfrak{Q}^{(2)}$ and $\mathfrak{P}^{(2)}$, we have that
\begin{enumerate}
\item[(i)]
$\mathfrak{P}^{(2)}$ is a head-constant echelonless second-order path that is not coherent and $i \in \bb{\mathfrak{P}^{(2)}}-1$ is the greatest index for which $\mathfrak{P}^{(2)0,i}$ is coherent.
\item[(ii)]
$\mathfrak{Q}^{(2)}$ is a head-constant echelonless second-order path.
\end{enumerate}

From (i) and taking into account Proposition~\ref{PDPthExt}, we have that the value of the path extension mapping of $\mathbf{f}^{(2)}$ at $\mathfrak{P}^{(2)}$ is given by
$$
f_{s}^{(2)\flat}\left(
\mathfrak{P}^{(2)}
\right)
=
f_{s}^{(2)\flat}\left(
\mathfrak{P}^{(2)i+1,\bb{\mathfrak{P}^{(2)}}-1}
\right)
\circ_{\varphi(s)}^{1\mathbf{Pth}_{\boldsymbol{\mathcal{B}}^{(2)}}}
f_{s}^{(2)\flat}\left(
\mathfrak{P}^{(2)0,i}
\right).
$$

Since $(\mathfrak{Q}^{(2)} \circ_{s}^{1\mathbf{Pth}_{\boldsymbol{\mathcal{A}}^{(2)}}} \mathfrak{P}^{(2)i+1,\bb{\mathfrak{P}^{(2)}}-1}, s)$$\prec_{\mathbf{Pth}_{\boldsymbol{\mathcal{A}}^{(2)}}}$$(\mathfrak{Q}^{(2)} \circ_{s}^{1\mathbf{Pth}_{\boldsymbol{\mathcal{A}}^{(2)}}} \mathfrak{P}^{(2)}, s)$ we have, by induction, that
\begin{align*}
&\left\llbracket
f_{s}^{(2)\flat}\left(
\mathfrak{Q}^{(2)}
\circ_{s}^{1\mathbf{Pth}_{\boldsymbol{\mathcal{A}}^{(2)}}}
\mathfrak{P}^{(2)i+1,\bb{\mathfrak{P}^{(2)}}-1}
\right)
\right\rrbracket_{\varphi(s)}
\\
&\hspace{3cm}
=
\left\llbracket
f_{s}^{(2)\flat}\left(
\mathfrak{Q}^{(2)}
\right)
\right\rrbracket_{\varphi(s)}
\circ_{\varphi(s)}^{1\mathbf{Pth}_{\boldsymbol{\mathcal{B}}^{(2)}}}
\left\llbracket
f_{s}^{(2)\flat}\left(
\mathfrak{P}^{(2)i+1,\bb{\mathfrak{P}^{(2)}}-1}
\right)
\right\rrbracket_{\varphi(s)}.
\end{align*}

So, considering the foregoing, we can affirm that
\begin{flushleft}
$
\left\llbracket
f_{s}^{(2)\flat}\left(
\mathfrak{Q}^{(2)} \circ_{s}^{1\mathbf{Pth}_{\boldsymbol{\mathcal{A}}^{(2)}}} \mathfrak{P}^{(2)}
\right)
\right\rrbracket_{\varphi(s)}
$
\allowdisplaybreaks
\begin{align*}
&=
\left\llbracket
f_{s}^{(2)\flat}\left(
\mathfrak{Q}^{(2)} \circ_{s}^{1\mathbf{Pth}_{\boldsymbol{\mathcal{A}}^{(2)}}} \mathfrak{P}^{(2)i+1,\bb{\mathfrak{P}^{(2)}}-1}
\right)
\circ_{\varphi(s)}^{1\mathbf{Pth}_{\boldsymbol{\mathcal{B}}^{(2)}}}
f_{s}^{(2)\flat}\left(
\mathfrak{P}^{(2)0,i}
\right)
\right\rrbracket_{\varphi(s)}
\tag{1}
\\
&=
\left\llbracket
f_{s}^{(2)\flat}\left(
\mathfrak{Q}^{(2)} \circ_{s}^{1\mathbf{Pth}_{\boldsymbol{\mathcal{A}}^{(2)}}} \mathfrak{P}^{(2)i+1,\bb{\mathfrak{P}^{(2)}}-1}
\right)
\circ_{s}^{1\mathbf{Pth}_{\boldsymbol{\mathcal{B}}^{(2)}}^{\mathbf{f}^{(2)}}}
f_{s}^{(2)\flat}\left(
\mathfrak{P}^{(2)0,i}
\right)
\right\rrbracket_{\varphi(s)}
\tag{2}
\\
&=
\left\llbracket
f_{s}^{(2)\flat}\left(
\mathfrak{Q}^{(2)} \circ_{s}^{1\mathbf{Pth}_{\boldsymbol{\mathcal{A}}^{(2)}}} \mathfrak{P}^{(2)i+1,\bb{\mathfrak{P}^{(2)}}-1}
\right)
\right\rrbracket_{\varphi(s)}
\circ_{s}^{1\llbracket\mathbf{Pth}_{\boldsymbol{\mathcal{B}}^{(2)}}^{\mathbf{f}^{(2)}}\rrbracket}
\\
&\hspace{8cm}
\left\llbracket
f_{s}^{(2)\flat}\left(
\mathfrak{P}^{(2)0,i}
\right)
\right\rrbracket_{\varphi(s)}
\tag{3}
\\
&=
\left(
\left\llbracket
f_{s}^{(2)\flat}\left(
\mathfrak{Q}^{(2)}
\right)
\right\rrbracket_{\varphi(s)}
\circ_{s}^{1\llbracket\mathbf{Pth}_{\boldsymbol{\mathcal{B}}^{(2)}}^{\mathbf{f}^{(2)}}\rrbracket} 
\left\llbracket
f_{s}^{(2)\flat}\left(
\mathfrak{P}^{(2)i+1,\bb{\mathfrak{P}}-1}
\right)
\right\rrbracket_{\varphi(s)}
\right)
\\
&\hspace{6.5cm}
\circ_{s}^{1\llbracket\mathbf{Pth}_{\boldsymbol{\mathcal{B}}^{(2)}}^{\mathbf{f}^{(2)}}\rrbracket}
\left\llbracket
f_{s}^{(2)\flat}\left(
\mathfrak{P}^{(2)0,i}
\right)
\right\rrbracket_{\varphi(s)}
\tag{4}
\\
&=
\left\llbracket
f_{s}^{(2)\flat}\left(
\mathfrak{Q}^{(2)}
\right)
\right\rrbracket_{\varphi(s)}
\circ_{s}^{1\llbracket\mathbf{Pth}_{\boldsymbol{\mathcal{B}}^{(2)}}^{\mathbf{f}^{(2)}}\rrbracket} 
\\
&\hspace{2cm}
\left(
\left\llbracket
f_{s}^{(2)\flat}\left(
\mathfrak{P}^{(2)i+1,\bb{\mathfrak{P}}-1}
\right)
\right\rrbracket_{\varphi(s)}
\circ_{s}^{1\llbracket\mathbf{Pth}_{\boldsymbol{\mathcal{B}}^{(2)}}^{\mathbf{f}^{(2)}}\rrbracket}
\left\llbracket
f_{s}^{(2)\flat}\left(
\mathfrak{P}^{(2)0,i}
\right)
\right\rrbracket_{\varphi(s)}
\right)
\tag{5}
\\
&=
\left\llbracket
f_{s}^{(2)\flat}\left(
\mathfrak{Q}^{(2)}
\right)
\right\rrbracket_{\varphi(s)}
\circ_{s}^{1\llbracket\mathbf{Pth}_{\boldsymbol{\mathcal{B}}^{(2)}}^{\mathbf{f}^{(2)}}\rrbracket} 
\\
&\hspace{3.5cm}
\left\llbracket
f_{s}^{(2)\flat}\left(
\mathfrak{P}^{(2)i+1,\bb{\mathfrak{P}^{(2)}}-1}
\right)
\circ_{s}^{1\mathbf{Pth}_{\boldsymbol{\mathcal{B}}^{(2)}}^{\mathbf{f}^{(2)}}}
f_{s}^{(2)\flat}\left(
\mathfrak{P}^{(2)0,i}
\right)
\right\rrbracket_{\varphi(s)}
\tag{6}
\\
&=
\left\llbracket
f_{s}^{(2)\flat}\left(
\mathfrak{Q}^{(2)}
\right)
\right\rrbracket_{\varphi(s)}
\circ_{s}^{1\llbracket\mathbf{Pth}_{\boldsymbol{\mathcal{B}}^{(2)}}^{\mathbf{f}^{(2)}}\rrbracket}
\\
&\hspace{3.5cm}
\left\llbracket
f_{s}^{(2)\flat}\left(
\mathfrak{P}^{(2)i+1,\bb{\mathfrak{P}}-1}
\right)
\circ_{\varphi(s)}^{1\mathbf{Pth}_{\boldsymbol{\mathcal{B}}^{(1)}}}
f_{s}^{(2)\flat}\left(
\mathfrak{P}^{(2)0,i}
\right)
\right\rrbracket_{\varphi(s)}
\tag{7}
\\
&=
\left\llbracket
f_{s}^{(2)\flat}\left(
\mathfrak{Q}^{(2)}
\right)
\right\rrbracket_{\varphi(s)}
\circ_{s}^{1\llbracket\mathbf{Pth}_{\boldsymbol{\mathcal{B}}^{(2)}}^{\mathbf{f}^{(2)}}\rrbracket}
\\
&\hspace{3.5cm}
\left\llbracket
f_{s}^{(2)\flat}\left(
\mathfrak{P}^{(2)i+1,\bb{\mathfrak{P}^{(2)}}-1}
\circ_{s}^{1\mathbf{Pth}_{\boldsymbol{\mathcal{A}}^{(2)}}}
\mathfrak{P}^{(2)0,i}
\right)
\right\rrbracket_{\varphi(s)}
\tag{8}
\\
&=
\left\llbracket
f_{s}^{(2)\flat}\left(
\mathfrak{Q}^{(2)}
\right)
\right\rrbracket_{\varphi(s)}
\circ_{s}^{1\llbracket\mathbf{Pth}_{\boldsymbol{\mathcal{B}}^{(2)}}^{\mathbf{f}^{(2)}}\rrbracket} 
\left\llbracket
f_{s}^{(2)\flat}\left(
\mathfrak{P}^{(2)}
\right)
\right\rrbracket_{\varphi(s)}.
\tag{9}
\end{align*}
\end{flushleft}

The first equality unravels the definition of the second-order path extension mapping $f^{(2)\flat}$ introduced in Proposition~\ref{PDPthExt};
the second equality recovers the definition of the operation symbol $\circ_{s}^{1}$ in the partial $\Sigma^{\boldsymbol{\mathcal{A}}^{(2)}}$-algebra $\mathbf{Pth}_{\boldsymbol{\mathcal{B}}^{(2)}}^{\mathbf{f}^{(2)}}$ introduced in Proposition~\ref{PDPthBCatAlg};
the third equality recovers the definition of the operation symbol $\circ_{s}^{1}$ in the partial $\Sigma^{\boldsymbol{\mathcal{A}}^{(1)}}$-algebra $\llbracket\mathbf{Pth}_{\boldsymbol{\mathcal{B}}^{(2)}}^{\mathbf{f}^{(2)}}\rrbracket$ introduced in Proposition~\ref{PDQPthBCatAlg};
the fourth equality follows by Artinian induction;
the fifth equality follows from the fact that, according to Proposition~\ref{PDPthComp}, the interpretation of the partial operation $\circ_{s}^{1}$ is associative;
the sixth equality unravels the definition of the operation symbol $\circ_{s}^{1}$ in the partial $\Sigma^{\boldsymbol{\mathcal{A}}^{(2)}}$-algebra $\llbracket\mathbf{Pth}_{\boldsymbol{\mathcal{B}}^{(2)}}^{\mathbf{f}^{(2)}}\rrbracket$ introduced in Proposition~\ref{PDQPthBCatAlg};
the seventh equality unravels the definition of the operation symbol $\circ_{s}^{1}$ in the partial $\Sigma^{\boldsymbol{\mathcal{A}}^{(2)}}$-algebra $\mathbf{Pth}_{\boldsymbol{\mathcal{B}}^{(2)}}^{\mathbf{f}^{(2)}}$ introduced in Proposition~\ref{PDPthBCatAlg};
the eight equality recovers the definition of the path extension mapping $f^{(2)\flat}$ introduced in Proposition~\ref{PPthExt};
finally, the last equality recovers the definition of $\mathfrak{P}^{(2)}$.

The case $i \in \bb{\mathfrak{P}^{(2)}}-1$ follows.

If $(2.2.2)$, i.e., if we find ourselves in the case where $\mathfrak{Q}^{(2)} \circ_{s}^{1\mathbf{Pth}_{\boldsymbol{\mathcal{A}}^{(2)}}} \mathfrak{P}^{(2)}$ is a head-constant echelonless second-order path that is not coherent and $i = \bb{\mathfrak{P}^{(2)}}-1$ is the greatest index for which $(\mathfrak{Q}^{(2)} \circ_{s}^{1\mathbf{Pth}_{\boldsymbol{\mathcal{A}}^{(2)}}} \mathfrak{P}^{(2)})^{0,i}$ is coherent then, regarding the second-order paths $\mathfrak{Q}^{(2)}$ and $\mathfrak{P}^{(2)}$, we have that
\begin{enumerate}
\item[(i)]
$\mathfrak{P}^{(2)}$ is a coherent head-constant echelonless second-order path.
\item[(ii)]
$\mathfrak{Q}^{(2)}$ is a head-constant echelonless second-order path.
\end{enumerate}

So, considering the foregoing, we can affirm that
\begin{flushleft}
$
\left\llbracket
f_{s}^{(2)\flat}\left(
\mathfrak{Q}^{(2)} \circ_{s}^{1\mathbf{Pth}_{\boldsymbol{\mathcal{A}}^{(2)}}} \mathfrak{P}^{(2)}
\right)
\right\rrbracket_{\varphi(s)}
$
\allowdisplaybreaks
\begin{align*}
&=
\left\llbracket
f_{s}^{(2)\flat}\left(
\mathfrak{Q}^{(2)}
\right)
\circ_{\varphi(s)}^{1\mathbf{Pth}_{\boldsymbol{\mathcal{B}}^{(2)}}}
f_{s}^{(2)\flat}\left(
\mathfrak{P}^{(2)}
\right)
\right\rrbracket_{\varphi(s)}
\tag{1}
\\
&=
\left\llbracket
f_{s}^{(2)\flat}\left(
\mathfrak{Q}^{(2)}
\right)
\circ_{s}^{1\mathbf{Pth}_{\boldsymbol{\mathcal{B}}^{(2)}}^{\mathbf{f}^{(2)}}}
f_{s}^{(2)\flat}\left(
\mathfrak{P}^{(2)}
\right)
\right\rrbracket_{\varphi(s)}
\tag{2}
\\
&=
\left\llbracket
f_{s}^{(2)\flat}\left(
\mathfrak{Q}^{(2)}
\right)
\right\rrbracket_{\varphi(s)}
\circ_{s}^{1\llbracket\mathbf{Pth}_{\boldsymbol{\mathcal{B}}^{(2)}}^{\mathbf{f}^{(2)}}\rrbracket}
\left\llbracket
f_{s}^{(2)\flat}\left(
\mathfrak{P}^{(2)}
\right)
\right\rrbracket_{\varphi(s)}
\tag{3}
\end{align*}
\end{flushleft}

The first equality unravels the definition of the second-order path extension mapping $f^{(2)\flat}$ introduced in Proposition~\ref{PDPthExt};
the second equality recovers the definition of the operation symbol $\circ_{s}^{1}$ in the partial $\Sigma^{\boldsymbol{\mathcal{A}}^{(2)}}$-algebra $\mathbf{Pth}_{\boldsymbol{\mathcal{B}}^{(2)}}^{\mathbf{f}^{(2)}}$ introduced in Proposition~\ref{PDPthBCatAlg};
finally, the last equality recovers the definition of the operation symbol $\circ_{s}^{1}$ in the partial $\Sigma^{\boldsymbol{\mathcal{A}}^{(1)}}$-algebra $\llbracket\mathbf{Pth}_{\boldsymbol{\mathcal{B}}^{(2)}}^{\mathbf{f}^{(2)}}\rrbracket$ introduced in Proposition~\ref{PDQPthBCatAlg}.

The case $i = \bb{\mathfrak{P}^{(2)}}-1$ follows.

If $(2.2.3)$, i.e., if we find ourselves in the case where $\mathfrak{Q}^{(2)} \circ_{s}^{1\mathbf{Pth}_{\boldsymbol{\mathcal{A}}^{(2)}}} \mathfrak{P}^{(2)}$ is a head-constant echelonless second-order path that is not coherent and $i \in [\bb{\mathfrak{P}^{(2)}}, \bb{\mathfrak{Q}^{(2)} \circ_{s}^{1\mathbf{Pth}_{\boldsymbol{\mathcal{A}}^{(2)}}} \mathfrak{P}^{(2)}}-1]$ is the greatest index for which $(\mathfrak{Q}^{(2)} \circ_{s}^{1\mathbf{Pth}_{\boldsymbol{\mathcal{A}}^{(2)}}} \mathfrak{P}^{(2)})^{0,i}$ is coherent then, regarding the second-order paths $\mathfrak{Q}^{(2)}$ and $\mathfrak{P}^{(2)}$, we have that
\begin{enumerate}
\item[(i)]
$\mathfrak{P}^{(2)}$ is a coherent head-constant echelonless second-order path.
\item[(ii)]
$\mathfrak{Q}^{(2)}$ is a head-constant echelonless second-order path that is not coherent and $j = i - \bb{\mathfrak{P}^{(2)}} \in \bb{\mathfrak{Q}^{(2)}}-1$ is the greatest index for which $\mathfrak{Q}^{(2)0,j}$ is coherent.
\end{enumerate}

From (ii) and taking into account Proposition~\ref{PDPthExt}, we have that the value of the path extension mapping of $\mathbf{f}^{(2)}$ at $\mathfrak{Q}^{(2)}$ is given by
$$
f_{s}^{(2)\flat}\left(
\mathfrak{Q}^{(2)}
\right)
=
f_{s}^{(2)\flat}\left(
\mathfrak{Q}^{(2)j+1,\bb{\mathfrak{Q}^{(2)}}-1}
\right)
\circ_{\varphi(s)}^{1\mathbf{Pth}_{\boldsymbol{\mathcal{B}}^{(2)}}}
f_{s}^{(2)\flat}\left(
\mathfrak{Q}^{(2)0,j}
\right).
$$

Since $(\mathfrak{Q}^{(2)0, j} \circ_{s}^{1\mathbf{Pth}_{\boldsymbol{\mathcal{A}}^{(2)}}} \mathfrak{P}^{(2)}, s)$$\prec_{\mathbf{Pth}_{\boldsymbol{\mathcal{A}}^{(2)}}}$$(\mathfrak{Q}^{(2)} \circ_{s}^{1\mathbf{Pth}_{\boldsymbol{\mathcal{A}}^{(2)}}} \mathfrak{P}^{(2)}, s)$ we have, by induction, that
\begin{align*}
&\left\llbracket
f_{s}^{(2)\flat}\left(
\mathfrak{Q}^{(2)0, j}
\circ_{s}^{1\mathbf{Pth}_{\boldsymbol{\mathcal{A}}^{(2)}}}
\mathfrak{P}^{(2)}
\right)
\right\rrbracket_{\varphi(s)}
\\
&\hspace{3cm}
=
\left\llbracket
f_{s}^{(2)\flat}\left(
\mathfrak{Q}^{(2)0, j}
\right)
\right\rrbracket_{\varphi(s)}
\circ_{s}^{1\llbracket\mathbf{Pth}_{\boldsymbol{\mathcal{B}}^{(2)}}^{\mathbf{f}^{(2)}}\rrbracket}
\left\llbracket
f_{s}^{(2)\flat}\left(
\mathfrak{P}^{(2)}
\right)
\right\rrbracket_{\varphi(s)}.
\end{align*}

So, considering the foregoing, we can affirm that
\begin{flushleft}
$
\left\llbracket
f_{s}^{(2)\flat}\left(
\mathfrak{Q}^{(2)} \circ_{s}^{1\mathbf{Pth}_{\boldsymbol{\mathcal{A}}^{(2)}}} \mathfrak{P}^{(2)}
\right)
\right\rrbracket_{\varphi(s)}
$
\allowdisplaybreaks
\begin{align*}
&=
\left\llbracket
f_{s}^{(2)\flat}\left(
\mathfrak{Q}^{(2)j+1, \bb{\mathfrak{Q}^{(2)}}-1}
\right)
\circ_{\varphi(s)}^{1\mathbf{Pth}_{\boldsymbol{\mathcal{B}}^{(2)}}}
f_{s}^{(2)\flat}\left(
\mathfrak{Q}^{(2)0, j}
\circ_{s}^{1\mathbf{Pth}_{\boldsymbol{\mathcal{A}}^{(2)}}}
\mathfrak{P}^{(2)}
\right)
\right\rrbracket_{\varphi(s)}
\tag{1}
\\
&=
\left\llbracket
f_{s}^{(2)\flat}\left(
\mathfrak{Q}^{(2)j+1, \bb{\mathfrak{Q}^{(2)}}-1}
\right)
\circ_{s}^{1\mathbf{Pth}_{\boldsymbol{\mathcal{B}}^{(2)}}^{\mathbf{f}^{(2)}}}
f_{s}^{(2)\flat}\left(
\mathfrak{Q}^{(2)0, j}
\circ_{s}^{1\mathbf{Pth}_{\boldsymbol{\mathcal{A}}^{(2)}}}
\mathfrak{P}^{(2)}
\right)
\right\rrbracket_{\varphi(s)}
\tag{2}
\\
&=
\left\llbracket
f_{s}^{(2)\flat}\left(
\mathfrak{Q}^{(2)j+1, \bb{\mathfrak{Q}^{(2)}}-1}
\right)
\right\rrbracket_{\varphi(s)}
\circ_{s}^{1\llbracket\mathbf{Pth}_{\boldsymbol{\mathcal{B}}^{(2)}}^{\mathbf{f}^{(2)}}\rrbracket}
\\
&\hspace{6cm}
\left\llbracket
f_{s}^{(2)\flat}\left(
\mathfrak{Q}^{(2)0, j}
\circ_{s}^{1\mathbf{Pth}_{\boldsymbol{\mathcal{A}}^{(2)}}}
\mathfrak{P}^{(2)}
\right)
\right\rrbracket_{\varphi(s)}
\tag{3}
\\
&=
\left\llbracket
f_{s}^{(2)\flat}\left(
\mathfrak{Q}^{(2)j+1, \bb{\mathfrak{Q}^{(2)}}-1}
\right)
\right\rrbracket_{\varphi(s)}
\circ_{s}^{1\llbracket\mathbf{Pth}_{\boldsymbol{\mathcal{B}}^{(2)}}^{\mathbf{f}^{(2)}}\rrbracket}
\\
&\hspace{3cm}
\left(
\left\llbracket
f_{s}^{(2)\flat}\left(
\mathfrak{Q}^{(2)0, j}
\right)
\right\rrbracket_{\varphi(s)}
\circ_{s}^{1\llbracket\mathbf{Pth}_{\boldsymbol{\mathcal{B}}^{(2)}}^{\mathbf{f}^{(2)}}\rrbracket}
\left\llbracket
f_{s}^{(2)\flat}\left(
\mathfrak{P}^{(2)}
\right)
\right\rrbracket_{\varphi(s)}
\right)
\tag{4}
\\
&=
\left(
\left\llbracket
f_{s}^{(2)\flat}\left(
\mathfrak{Q}^{(2)j+1, \bb{\mathfrak{Q}^{(2)}}-1}
\right)
\right\rrbracket_{\varphi(s)}
\circ_{s}^{1\llbracket\mathbf{Pth}_{\boldsymbol{\mathcal{B}}^{(2)}}^{\mathbf{f}^{(2)}}\rrbracket}
\left\llbracket
f_{s}^{(2)\flat}\left(
\mathfrak{Q}^{(2)0, j}
\right)
\right\rrbracket_{\varphi(s)}
\right)
\\
&\hspace{7cm}
\circ_{s}^{1\llbracket\mathbf{Pth}_{\boldsymbol{\mathcal{B}}^{(2)}}^{\mathbf{f}^{(2)}}\rrbracket}
\left\llbracket
f_{s}^{(2)\flat}\left(
\mathfrak{P}^{(2)}
\right)
\right\rrbracket_{\varphi(s)}
\tag{5}
\\
&=
\left\llbracket
f_{s}^{(2)\flat}\left(
\mathfrak{Q}^{(2)j+1, \bb{\mathfrak{Q}^{(2)}}-1}
\right)
\circ_{s}^{1\mathbf{Pth}_{\boldsymbol{\mathcal{B}}^{(2)}}^{\mathbf{f}^{(2)}}}
f_{s}^{(2)\flat}\left(
\mathfrak{Q}^{(2)0, j}
\right)
\right\rrbracket_{\varphi(s)}
\\
&\hspace{7cm}
\circ_{s}^{1\llbracket\mathbf{Pth}_{\boldsymbol{\mathcal{B}}^{(2)}}^{\mathbf{f}^{(2)}}\rrbracket}
\left\llbracket
f_{s}^{(2)\flat}\left(
\mathfrak{P}^{(2)}
\right)
\right\rrbracket_{\varphi(s)}
\tag{6}
\\
&=
\left\llbracket
f_{s}^{(2)\flat}\left(
\mathfrak{Q}^{(2)j+1, \bb{\mathfrak{Q}^{(2)}}-1}
\right)
\circ_{\varphi(s)}^{1\mathbf{Pth}_{\boldsymbol{\mathcal{B}}^{(2)}}}
f_{s}^{(2)\flat}\left(
\mathfrak{Q}^{(2)0, j}
\right)
\right\rrbracket_{\varphi(s)}
\\
&\hspace{7cm}
\circ_{s}^{1\llbracket\mathbf{Pth}_{\boldsymbol{\mathcal{B}}^{(2)}}^{\mathbf{f}^{(2)}}\rrbracket}
\left\llbracket
f_{s}^{(2)\flat}\left(
\mathfrak{P}^{(2)}
\right)
\right\rrbracket_{\varphi(s)}
\tag{7}
\\
&=
\left\llbracket
f_{s}^{(2)\flat}\left(
\mathfrak{Q}^{(2)j+1, \bb{\mathfrak{Q}^{(2)}}-1}
\circ_{s}^{1\mathbf{Pth}_{\boldsymbol{\mathcal{A}}^{(2)}}}
\mathfrak{Q}^{(2)0, j}
\right)
\right\rrbracket_{\varphi(s)}
\\
&\hspace{7cm}
\circ_{s}^{1\llbracket\mathbf{Pth}_{\boldsymbol{\mathcal{B}}^{(2)}}^{\mathbf{f}^{(2)}}\rrbracket}
\left\llbracket
f_{s}^{(2)\flat}\left(
\mathfrak{P}^{(2)}
\right)
\right\rrbracket_{\varphi(s)}
\tag{8}
\\
&=
\left\llbracket
f_{s}^{(2)\flat}\left(
\mathfrak{Q}^{(2)}
\right)
\right\rrbracket_{\varphi(s)}
\circ_{s}^{1\llbracket\mathbf{Pth}_{\boldsymbol{\mathcal{B}}^{(2)}}^{\mathbf{f}^{(2)}}\rrbracket} 
\left\llbracket
f_{s}^{(2)\flat}\left(
\mathfrak{P}^{(2)}
\right)
\right\rrbracket_{\varphi(s)}.
\tag{9}
\end{align*}
\end{flushleft}

The first equality unravels the definition of the second-order path extension mapping $f^{(2)\flat}$ introduced in Proposition~\ref{PDPthExt};
the second equality recovers the definition of the operation symbol $\circ_{s}^{1}$ in the partial $\Sigma^{\boldsymbol{\mathcal{A}}^{(2)}}$-algebra $\mathbf{Pth}_{\boldsymbol{\mathcal{B}}^{(2)}}^{\mathbf{f}^{(2)}}$ introduced in Proposition~\ref{PDPthBCatAlg};
the third equality recovers the definition of the operation symbol $\circ_{s}^{1}$ in the partial $\Sigma^{\boldsymbol{\mathcal{A}}^{(1)}}$-algebra $\llbracket\mathbf{Pth}_{\boldsymbol{\mathcal{B}}^{(2)}}^{\mathbf{f}^{(2)}}\rrbracket$ introduced in Proposition~\ref{PDQPthBCatAlg};
the fourth equality follows by Artinian induction;
the fifth equality follows from the fact that, according to Proposition~\ref{PDPthComp}, the interpretation of the partial operation $\circ_{s}^{1}$ is associative;
the sixth equality unravels the definition of the operation symbol $\circ_{s}^{1}$ in the partial $\Sigma^{\boldsymbol{\mathcal{A}}^{(2)}}$-algebra $\llbracket\mathbf{Pth}_{\boldsymbol{\mathcal{B}}^{(2)}}^{\mathbf{f}^{(2)}}\rrbracket$ introduced in Proposition~\ref{PDQPthBCatAlg};
the seventh equality unravels the definition of the operation symbol $\circ_{s}^{1}$ in the partial $\Sigma^{\boldsymbol{\mathcal{A}}^{(2)}}$-algebra $\mathbf{Pth}_{\boldsymbol{\mathcal{B}}^{(2)}}^{\mathbf{f}^{(2)}}$ introduced in Proposition~\ref{PDPthBCatAlg};
the eight equality recovers the definition of the path extension mapping $f^{(2)\flat}$ introduced in Proposition~\ref{PPthExt};
finally, the last equality recovers the definition of $\mathfrak{Q}^{(2)}$.

The case $i \in [\bb{\mathfrak{P}^{(2)}}, \bb{\mathfrak{Q}^{(2)} \circ_{s}^{1\mathbf{Pth}_{\boldsymbol{\mathcal{A}}^{(2)}}} \mathfrak{P}^{(2)}}-1]$ follows.

This completes the case $(2.2)$.

If $(2.3)$, i.e., if $\mathfrak{Q}^{(2)} \circ_{s}^{1\mathbf{Pth}_{\boldsymbol{\mathcal{A}}^{(2)}}} \mathfrak{P}^{(2)}$ is a coherent head-constant echelonless second-order path then, regarding the paths $\mathfrak{Q}^{(2)}$ and $\mathfrak{P}^{(2)}$, we have that
\begin{enumerate}
\item[(i)]
$\mathfrak{P}^{(2)}$ is a coherent head-constant echelonless second-order path.
\item[(ii)]
$\mathfrak{Q}^{(2)}$ is a coherent head-constant echelonless second-order path.
\end{enumerate}

Therefore, for a unique word $\mathbf{s} \in S^{\star} - \{\lambda\}$ and a unique operation symbol $\tau \in \Sigma^{\boldsymbol{\mathcal{A}}^{(1)}}_{\mathbf{s}, s}$, the family of first-order translations occurring in $\mathfrak{P}^{(2)}$ is a family of fist-order translations of type $\tau$.

Let $(\mathfrak{P}_{j}^{(2)})_{j\in\bb{\mathbf{s}}}$ be the family of second-order paths we can extract from $\mathfrak{P}^{(2)}$ in virtue of Lemma~\ref{LDPthExtract}. Then, according to Proposition~\ref{PDPthExt}, we have that the value of the path extension mapping of $\mathbf{f}^{(2)}$ at $\mathfrak{P}^{(2)}$ is given by
$$
f_{s}^{(2)\flat} \left(\mathfrak{P}^{(2)}\right)
=
\tau^{\mathbf{Pth}_{\boldsymbol{\mathcal{B}}^{(2)}}^{\mathbf{f}^{(2)}}}\left(\left(
f_{s_{j}}^{(2)\flat}\left(
\mathfrak{P}_{j}^{(2)}
\right)
\right)_{j\in\bb{\mathbf{s}}}\right).
$$

Moreover, for the unique word $\mathbf{s} \in S^{\star} - \{\lambda\}$ and the unique operation symbol $\tau \in \Sigma^{\boldsymbol{\mathcal{A}}^{(1)}}_{\mathbf{s}, s}$, the family of first-order translations occurring in $\mathfrak{Q}^{(2)}$ is a family of fist-order translations of type $\tau$. Note that the operation symbol $\tau$ is the same as in case (i), since $\mathfrak{Q}^{(2)} \circ_{s}^{1\mathbf{Pth}_{\boldsymbol{\mathcal{A}}^{(2)}}} \mathfrak{P}^{(2)}$ is head-constant by hypothesis.

Let $(\mathfrak{Q}_{j}^{(2)})_{j\in\bb{\mathbf{s}}}$ be the family of second-order paths we can extract from $\mathfrak{Q}^{(2)}$ in virtue of Lemma~\ref{LDPthExtract}. Then, according to Proposition~\ref{PDPthExt}, we have that the value of the path extension mapping of $\mathbf{f}^{(2)}$ at $\mathfrak{Q}^{(2)}$ is given by
$$
f_{s}^{(2)\flat} \left(\mathfrak{Q}^{(2)}\right)
=
\tau^{\mathbf{Pth}_{\boldsymbol{\mathcal{B}}^{(2)}}^{\mathbf{f}^{(2)}}}\left(\left(
f_{s_{j}}^{(2)\flat}\left(
\mathfrak{Q}_{j}^{(2)}
\right)
\right)_{j\in\bb{\mathbf{s}}}\right).
$$

Let us consider $((\mathfrak{Q}^{(2)} \circ_{s}^{1\mathbf{Pth}_{\boldsymbol{\mathcal{A}}^{(2)}}} \mathfrak{P}^{(2)})_{j})_{j\in\bb{\mathbf{s}}}$ the family of paths we can extract, in virtue of Lemma~\ref{LDPthExtract}, from $\mathfrak{Q}^{(2)} \circ_{s}^{1\mathbf{Pth}_{\boldsymbol{\mathcal{A}}^{(2)}}} \mathfrak{P}^{(2)}$. Let us note that, for every $j\in\bb{\mathbf{s}}$, it is the case that
$$
(\mathfrak{Q}^{(2)} \circ_{s}^{1\mathbf{Pth}_{\boldsymbol{\mathcal{A}}^{(1)}}} \mathfrak{P}^{(2)})_{j}
=
\mathfrak{Q}_{j}^{(2)} \circ_{s_{j}}^{1\mathbf{Pth}_{\boldsymbol{\mathcal{A}}^{(2)}}} \mathfrak{P}_{j}^{(2)}.
$$
Then according to Proposition~\ref{PDPthExt}, we have that the value of the path extension mapping of $\mathbf{f}^{(2)}$ at $\mathfrak{Q}^{(2)} \circ_{s}^{1\mathbf{Pth}_{\boldsymbol{\mathcal{A}}^{(2)}}} \mathfrak{P}^{(2)}$ is given by
$$
f_{s}^{(2)\flat} \left(
\mathfrak{Q}^{(2)}
\circ_{s}^{1\mathbf{Pth}_{\boldsymbol{\mathcal{A}}^{(2)}}}
\mathfrak{P}^{(2)}
\right)
=
\tau^{\mathbf{Pth}_{\boldsymbol{\mathcal{B}}^{(2)}}^{\mathbf{f}^{(2)}}}\left(\left(
f_{s_{j}}^{(2)\flat}\left(
\mathfrak{Q}_{j}^{(2)} \circ_{s_{j}}^{1\mathbf{Pth}_{\boldsymbol{\mathcal{A}}^{(2)}}} \mathfrak{P}_{j}^{(2)}
\right)
\right)_{j\in\bb{\mathbf{s}}}\right).
$$

Finally, note that, for every $j\in\bb{\mathbf{s}}$, $(\mathfrak{Q}_{j}^{(2)} \circ_{s}^{1\mathbf{Pth}_{\boldsymbol{\mathcal{A}}^{(2)}}} \mathfrak{P}_{j}^{(2)}, s_{j}) \prec_{\mathbf{Pth}_{\boldsymbol{\mathcal{A}}^{(2)}}} (\mathfrak{Q}^{(2)} \circ_{s}^{1\mathbf{Pth}_{\boldsymbol{\mathcal{A}}^{(2)}}} \mathfrak{P}^{(2)}, s)$. Thus, by induction, we have that
$$
\left\llbracket
f_{s_{j}}^{(2)\flat}\left(
\mathfrak{Q}_{j}^{(2)} \circ_{s_{j}}^{1\mathbf{Pth}_{\boldsymbol{\mathcal{A}}^{(2)}}} \mathfrak{P}_{j}^{(2)}
\right)
\right\rrbracket_{\varphi(s_{j})}
=
\left\llbracket
f_{s_{j}}^{(2)\flat}\left(
\mathfrak{Q}_{j}^{(2)}
\right)
\right\rrbracket_{\varphi(s_{j})}
\circ_{s_{j}}^{1\llbracket\mathbf{Pth}_{\boldsymbol{\mathcal{B}}^{(2)}}^{\mathbf{f}^{(2)}}\rrbracket}
\left\llbracket
f_{s_{j}}^{(2)\flat}\left(
\mathfrak{P}_{j}^{(2)}
\right)
\right\rrbracket_{\varphi(s_{j})}.
$$

The following chain of equalities holds
\begin{flushleft}
$
\left\llbracket
f_{s}^{(2)\flat}\left(
\mathfrak{Q}^{(2)} \circ_{s}^{1\mathbf{Pth}_{\boldsymbol{\mathcal{A}}^{(2)}}} \mathfrak{P}^{(2)}
\right)
\right\rrbracket_{\varphi(s)}
$
\allowdisplaybreaks
\begin{align*}
&=
\left\llbracket
\tau^{\mathbf{Pth}_{\boldsymbol{\mathcal{B}}^{(2)}}^{\mathbf{f}^{(2)}}}\left(\left(
f_{s_{j}}^{(2)\flat}\left(
\mathfrak{Q}_{j}^{(2)} \circ_{s_{j}}^{1\mathbf{Pth}_{\boldsymbol{\mathcal{A}}^{(2)}}} \mathfrak{P}_{j}^{(2)}
\right)
\right)_{j\in\bb{\mathbf{s}}}\right)
\right\rrbracket_{\varphi(s)}
\tag{1}
\\
&=
\tau^{\llbracket\mathbf{Pth}_{\boldsymbol{\mathcal{B}}^{(2)}}^{\mathbf{f}^{(2)}}\rrbracket} \left(\left(
\left\llbracket
f_{s_{j}}^{(2)\flat}\left(
\mathfrak{Q}_{j}^{(2)} \circ_{s_{j}}^{1\mathbf{Pth}_{\boldsymbol{\mathcal{A}}^{(2)}}} \mathfrak{P}_{j}^{(2)}
\right)
\right\rrbracket_{\varphi(s_{j})}
\right)_{j\in\bb{\mathbf{s}}}\right)
\tag{2}
\\
&=
\resizebox{.89\textwidth}{!}{%
$
\tau^{\llbracket\mathbf{Pth}_{\boldsymbol{\mathcal{B}}^{(2)}}^{\mathbf{f}^{(2)}}\rrbracket} \left(\left(
\left\llbracket
f_{s_{j}}^{(2)\flat}\left(
\mathfrak{Q}_{j}^{(2)}
\right)
\right\rrbracket_{\varphi(s_{j})}
\circ_{s_{j}}^{1\llbracket\mathbf{Pth}_{\boldsymbol{\mathcal{B}}^{(2)}}^{\mathbf{f}^{(2)}}\rrbracket}
\left\llbracket
f_{s_{j}}^{(2)\flat}\left(
\mathfrak{P}_{j}^{(2)}
\right)
\right\rrbracket_{\varphi(s_{j})}
\right)_{j\in\bb{\mathbf{s}}}\right)
$
}
\tag{3}
\\
&=
\tau^{\llbracket\mathbf{Pth}_{\boldsymbol{\mathcal{B}}^{(2)}}^{\mathbf{f}^{(2)}}\rrbracket} \left(\left(
\left\llbracket
f_{s_{j}}^{(2)\flat}\left(
\mathfrak{Q}_{j}^{(2)}
\right)
\right\rrbracket_{\varphi(s_{j})}
\right)_{j\in\bb{\mathbf{s}}}\right)
\circ_{s}^{1\llbracket\mathbf{Pth}_{\boldsymbol{\mathcal{B}}^{(2)}}^{\mathbf{f}^{(2)}}\rrbracket}
\\
&\hspace{5cm}
\tau^{\llbracket\mathbf{Pth}_{\boldsymbol{\mathcal{B}}^{(2)}}^{\mathbf{f}^{(2)}}\rrbracket} \left(\left(
\left\llbracket
f_{s_{j}}^{(2)\flat}\left(
\mathfrak{P}_{j}^{(2)}
\right)
\right\rrbracket_{\varphi(s_{j})}
\right)_{j\in\bb{\mathbf{s}}}\right)
\tag{4}
\\
&=
\left\llbracket
\tau^{\mathbf{Pth}_{\boldsymbol{\mathcal{B}}^{(2)}}^{\mathbf{f}^{(2)}}}\left(\left(
f_{s_{j}}^{(2)\flat}\left(
\mathfrak{Q}_{j}^{(2)}
\right)
\right)_{j\in\bb{\mathbf{s}}}\right)
\right\rrbracket_{\varphi(s)}
\circ_{s}^{1\llbracket\mathbf{Pth}_{\boldsymbol{\mathcal{B}}^{(2)}}^{\mathbf{f}^{(2)}}\rrbracket}
\\
&\hspace{5.5cm}
\left\llbracket
\tau^{\mathbf{Pth}_{\boldsymbol{\mathcal{B}}^{(2)}}^{\mathbf{f}^{(2)}}}\left(\left(
f_{s_{j}}^{(2)\flat}\left(
\mathfrak{P}_{j}^{(2)}
\right)
\right)_{j\in\bb{\mathbf{s}}}\right)
\right\rrbracket_{\varphi(s)}
\tag{5}
\\
&=
\left\llbracket
f_{s}^{(2)\flat}\left(
\mathfrak{Q}^{(2)}
\right)
\right\rrbracket_{\varphi(s)}
\circ_{s}^{0[\mathbf{Pth}_{\boldsymbol{\mathcal{B}}^{(1)}}^{\mathbf{f}^{(1)}}]}
\left\llbracket
f_{s}^{(1)\flat}\left(
\mathfrak{P}^{(2)}
\right)
\right\rrbracket_{\varphi(s)}.
\tag{6}
\end{align*}
\end{flushleft}

The first equality unravels the definition of the second-order path extension mapping $f^{(2)\flat}$ introduced in Proposition~\ref{PDPthExt} at a coherent head-constant echelonless second-order path;
the second equality recovers the interpretation of the operation symbol $\tau$ in the partial $\Sigma^{\boldsymbol{\mathcal{A}}^{(1)}}$-algebra $\llbracket\mathbf{Pth}_{\boldsymbol{\mathcal{B}}^{(2)}}^{\mathbf{f}^{(2)}(1,2)}\rrbracket$;
the third equality follows by Artinian induction;
the fourth equality follows from Proposition~\ref{PDQPthBVarB8} if $\tau \in \Sigma_{\mathbf{s}, s}$ or from Axiom~AB3 if $\tau = \circ_{s}^{0}$ since, by Proposition~\ref{PDPthVar} $\llbracket\mathbf{Pth}_{\boldsymbol{\mathcal{B}}^{(2)}}\rrbracket$ is a many-sorted partial $\Lambda^{\boldsymbol{\mathcal{B}}^{(2)}}$-algebra in $\mathbf{PAlg}(\boldsymbol{\mathcal{E}}^{\boldsymbol{\mathcal{B}}^{(2)}})$;
the fifth equality unravels the interpretation of the operation symbol $\tau$ in the partial $\Sigma^{\boldsymbol{\mathcal{A}}^{(1)}}$-algebra $\llbracket\mathbf{Pth}_{\boldsymbol{\mathcal{B}}^{(2)}}^{\mathbf{f}^{(2)}(1,2)}\rrbracket$ introduced in Proposition~\ref{PDQPthBCatAlg};
finally, the last equality recovers the definition of $f^{(2)\flat}$ at a coherent head-constant echelonless second-order path.

This completes case $(2.3)$.

This completes case $(2)$.

Hence, $\mathrm{pr}_{\boldsymbol{\mathcal{B}}^{(2)}, \varphi}^{\llbracket\cdot\rrbracket} \circ f^{(2)\flat}$ is compatible with the $1$-composition operation.

All in all, we conclude that $\mathrm{pr}_{\boldsymbol{\mathcal{B}}^{(2)}, \varphi}^{\llbracket\cdot\rrbracket} \circ f^{(2)\flat}$ is a $\Sigma^{\boldsymbol{\mathcal{A}}^{(2)}}$-homomorphism.

It remains to prove that
$$
\mathrm{Ker}(\mathrm{CH}^{(2)})
\vee
\Upsilon^{[1]}
\subseteq
\mathrm{Ker}(\mathrm{pr}_{\boldsymbol{\mathcal{B}}^{(2)}, \varphi}^{\llbracket\cdot\rrbracket} \circ f^{(2)\flat}).
$$

In order to prove that $\mathrm{Ker}(\mathrm{CH}^{(2)}_{\boldsymbol{\mathcal{A}}^{(2)}}) \vee \Upsilon^{[1]}$ is included in $\mathrm{Ker}(\mathrm{pr}_{\boldsymbol{\mathcal{B}}^{(2)}, \varphi}^{\llbracket\cdot\rrbracket} \circ f^{(2)\flat})$, it suffices to prove that both $\mathrm{Ker}(\mathrm{CH}^{(2)}_{\boldsymbol{\mathcal{A}}^{(2)}})$ and $\Upsilon^{[1]}$ are included in $\mathrm{Ker}(\mathrm{pr}_{\boldsymbol{\mathcal{B}}^{(2)}, \varphi}^{\llbracket\cdot\rrbracket} \circ f^{(2)\flat})$.

{\sffamily (i) $\mathrm{Ker}(\mathrm{CH}^{(2)}_{\boldsymbol{\mathcal{A}}^{(2)}}) \subseteq \mathrm{Ker}(\mathrm{pr}_{\boldsymbol{\mathcal{B}}^{(2)}, \varphi}^{\llbracket\cdot\rrbracket} \circ f^{(2)\flat})$.}

Let $s$ be a sort in $S$ and let $\mathfrak{Q}^{(2)}, \mathfrak{P}^{(2)}$ be two second-order paths in $\mathrm{Pth}_{\boldsymbol{\mathcal{A}}^{(2)}, s}$ satisfying that $(\mathfrak{Q}^{(2)}, \mathfrak{P}^{(2)}) \in \mathrm{Ker}(\mathrm{CH}^{(2)}_{\boldsymbol{\mathcal{A}}^{(2)}})_{s}$.

we will prove that
$$
\left\llbracket
f^{(2)\flat}_{s}\left(
\mathfrak{P}^{(2)}
\right)
\right\rrbracket_{\varphi(s)}
=
\left\llbracket
f^{(2)\flat}_{s}\left(
\mathfrak{Q}^{(2)}
\right)
\right\rrbracket_{\varphi(s)}
$$
by Artinian induction on $(\coprod \mathrm{Pth}_{\boldsymbol{\mathcal{A}}^{(2)}}, \leq_{\mathbf{Pth}_{\boldsymbol{\mathcal{A}}^{(2)}}})$.

{\sffamily Base step of the Artinian induction}

Let $(\mathfrak{P}^{(2)}, s)$ be a minimal element of $(\coprod \mathrm{Pth}_{\boldsymbol{\mathcal{A}}^{(2)}}, \leq_{\mathbf{Pth}_{\boldsymbol{\mathcal{A}}^{(2)}}})$. Then, by Proposition~\ref{PDMinimal}, the path $\mathfrak{P}^{(2)}$ is either (1) a $(2,[1])$-identity second-order path or (2) a second-order echelon.

In either case, the equality
$$
\left\llbracket
f^{(2)\flat}_{s}\left(
\mathfrak{P}^{(2)}
\right)
\right\rrbracket_{\varphi(s)}
=
\left\llbracket
f^{(2)\flat}_{s}\left(
\mathfrak{Q}^{(2)}
\right)
\right\rrbracket_{\varphi(s)}
$$
follows  directly from the fact that, according to Corollary~\ref{CDCHUId} or Proposition~\ref{PDCHEch}, $\mathfrak{P}^{(2)} = \mathfrak{Q}^{(2)}$.

{\sffamily Inductive step of the Artinian induction}

Let $(\mathfrak{P}^{(2)}, s)$ be a non-minimal element of $(\coprod\mathrm{Pth}_{\boldsymbol{\mathcal{A}}^{(2)}}, \leq_{\mathbf{Pth}_{\boldsymbol{\mathcal{A}}^{(2)}}})$. Let us suppose that, for every sort $t \in S$ and every second-order path $\mathfrak{P}'^{(2)}$ in $\mathrm{Pth}_{\boldsymbol{\mathcal{A}}^{(2)}, t}$, if $(\mathfrak{P}'^{(2)},t) <_{\mathbf{Pth}_{\boldsymbol{\mathcal{A}}^{(2)}}} (\mathfrak{P}^{(2)}, s)$, the the statement holds for $\mathfrak{P}'^{(2)}$, i.e., for every second-order path $\mathfrak{Q}'^{(2)}$ in $\mathrm{Pth}_{\boldsymbol{\mathcal{A}}^{(2)}, t}$, if $(\mathfrak{P}'^{(2)}, \mathfrak{Q}'^{(2)})\in\mathrm{Ker}(\mathrm{CH}^{(2)})_{t}$, then 
$$
\left\llbracket
f^{(2)\flat}_{t}\left(
\mathfrak{P}'^{(2)}
\right)
\right\rrbracket_{\varphi(t)}
=
\left\llbracket
f^{(2)\flat}_{t}\left(
\mathfrak{Q}'^{(2)}
\right)
\right\rrbracket_{\varphi(t)}.
$$

Let $(\mathfrak{P}^{(2)}, s)$ be a non-minimal element of $(\coprod \mathrm{Pth}_{\boldsymbol{\mathcal{A}}^{(2)}, \leq_{\mathbf{Pth}_{\boldsymbol{\mathcal{A}}^{(2)}}}})$. We can assume that $\mathfrak{P}^{(2)}$ is not a $(2,[1])$-identity second-order path, since that case has already been considered. By Lemma~\ref{LDOrdI}, that $\mathfrak{P}^{(2)}$ is either (1) a second-order path of length strictly greater than one containing at least one second-order echelon or (2) an echelonless second-order path.

If (1), then let $i \in \bb{\mathfrak{P}^{(2)}}$ be the first index for which the one-step subpath $\mathfrak{P}^{(2), i,i}$ is a second-order echelon. We distinguish two cases accordingly.

If $i=0$, i.e., if $\mathfrak{P}^{(2)}$ has its first second-order echelon on its first step, then according to Proposition~\ref{PDPthExt}, we have that
$$
f^{(2)\flat}_{s}\left(
\mathfrak{P}^{(2)}
\right)
=
f^{(2)\flat}_{s}\left(
\mathfrak{P}^{(2),1,\bb{\mathfrak{P}^{(2)}}-1}
\right)
\circ_{\varphi(s)}^{1\mathbf{Pth}_{\boldsymbol{\mathcal{B}}^{(2)}}}
f^{(2)\flat}_{s}\left(
\mathfrak{P}^{(2),0,0}
\right).
$$

Since $\mathrm{CH}^{(2)}_{\boldsymbol{\mathcal{A}}^{(2)}, s}(\mathfrak{P}^{(2)}) \in \eta^{(2,\mathcal{A}^{(2)})}[\mathcal{A}^{(2)}]_{s}^{\mathrm{int}}$ and $(\mathfrak{P}^{(2)}, \mathfrak{Q}^{(2)}) \in \mathrm{Ker}(\mathrm{CH}^{(2)}_{\boldsymbol{\mathcal{A}}^{(2)}})_{s}$, we have, by Lemma~\ref{LDCHEchInt}, that $\mathfrak{Q}^{(2)}$ is a path of length strictly grater than one containing its first second-order echelon on its first step.

Thus, according to Proposition~\ref{PDPthExt}, we have that
$$
f^{(2)\flat}_{s}\left(
\mathfrak{Q}^{(2)}
\right)
=
f^{(2)\flat}_{s}\left(
\mathfrak{Q}^{(2),1,\bb{\mathfrak{P}^{(2)}}-1}
\right)
\circ_{\varphi(s)}^{1\mathbf{Pth}_{\boldsymbol{\mathcal{B}}^{(2)}}}
f^{(2)\flat}_{s}\left(
\mathfrak{Q}^{(2),0,0}
\right).
$$

Since $(\mathfrak{P}^{(2)}, \mathfrak{Q}^{(2)}) \in \mathrm{Ker}(\mathrm{CH}^{(2)}_{\boldsymbol{\mathcal{A}}^{(2)}})_{s}$, we have that $(\mathfrak{P}^{(2),1,\bb{\mathfrak{P}^{(2)}}-1}, \mathfrak{Q}^{(2),1,\bb{\mathfrak{Q}^{(2)}}-1}$ and $(\mathfrak{P}^{(2),0,0}, \mathfrak{Q}^{(2),0,0})$ are in $\mathrm{Ker}(\mathrm{CH}^{(2)})_{s}$. Note that, according to Definition~\ref{DDOrd}, we have that $(\mathfrak{P}^{(2),0,0}, s)$ and $(\mathfrak{P}^{(2),1,\bb{\mathfrak{P}^{(2)}}-1}, s)$$\prec_{\mathbf{Pth}_{\boldsymbol{\mathcal{A}}^{(2)}}}$-precede $(\mathfrak{P}^{(2)}, s)$.

Therefore, by the inductive hypothesis, the second-order paths $\mathfrak{P}^{(2),0,0}$ and $\mathfrak{Q}^{(2),0,0}$, and the second-order paths $\mathfrak{P}^{(2),1,\bb{\mathfrak{P}^{(2)}}-1}$ and $\mathfrak{Q}^{(2),1,\bb{\mathfrak{Q}^{(2)}}-1}$ satisfy
\begin{align*}
\left\llbracket
f^{(2)\flat}_{s}\left(
\mathfrak{P}^{(2),0,0}
\right)
\right\rrbracket_{\varphi(s)}
&=
\left\llbracket
f^{(2)\flat}_{s}\left(
\mathfrak{Q}^{(2),0,0}
\right)
\right\rrbracket_{\varphi(s)}
\\
\left\llbracket
f^{(2)\flat}_{s}\left(
\mathfrak{P}^{(2),1,\bb{\mathfrak{P}^{(2)}}-1}
\right)
\right\rrbracket_{\varphi(s)}
&=
\left\llbracket
f^{(2)\flat}_{s}\left(
\mathfrak{Q}^{(2),1,\bb{\mathfrak{Q}^{(2)}}-1}
\right)
\right\rrbracket_{\varphi(s)}.
\end{align*}

Thus, the following chain of equalities holds
\begin{flushleft}
$
\left\llbracket
f^{(2)\flat}_{s}\left(
\mathfrak{P}^{(2)}
\right)
\right\rrbracket_{\varphi(s)}
$
\allowdisplaybreaks
\begin{align*}
&=
\left\llbracket
f^{(2)\flat}_{s}\left(
\mathfrak{P}^{(2),1,\bb{\mathfrak{P}^{(2)}}-1}
\right)
\circ_{\varphi(s)}^{1\mathbf{Pth}_{\boldsymbol{\mathcal{B}}^{(2)}}}
f^{(2)\flat}_{s}\left(
\mathfrak{P}^{(2),0,0}
\right)
\right\rrbracket_{\varphi(s)}
\tag{1}
\\
&=
\left\llbracket
f^{(2)\flat}_{s}\left(
\mathfrak{P}^{(2),1,\bb{\mathfrak{P}^{(2)}}-1}
\right)
\right\rrbracket_{\varphi(s)}
\circ_{\varphi(s)}^{1\llbracket\mathbf{Pth}_{\boldsymbol{\mathcal{B}}^{(2)}}\rrbracket}
\left\llbracket
f^{(2)\flat}_{s}\left(
\mathfrak{P}^{(2),0,0}
\right)
\right\rrbracket_{\varphi(s)}
\tag{2}
\\
&=
\left\llbracket
f^{(2)\flat}_{s}\left(
\mathfrak{Q}^{(2),1,\bb{\mathfrak{Q}^{(2)}}-1}
\right)
\right\rrbracket_{\varphi(s)}
\circ_{\varphi(s)}^{1\llbracket\mathbf{Pth}_{\boldsymbol{\mathcal{B}}^{(2)}}\rrbracket}
\left\llbracket
f^{(2)\flat}_{s}\left(
\mathfrak{Q}^{(2),0,0}
\right)
\right\rrbracket_{\varphi(s)}
\tag{3}
\\
&=
\left\llbracket
f^{(2)\flat}_{s}\left(
\mathfrak{Q}^{(2),1,\bb{\mathfrak{Q}^{(2)}}-1}
\right)
\circ_{\varphi(s)}^{1\mathbf{Pth}_{\boldsymbol{\mathcal{B}}^{(2)}}}
f^{(2)\flat}_{s}\left(
\mathfrak{Q}^{(2),0,0}
\right)
\right\rrbracket_{\varphi(s)}
\tag{4}
\\
&=
\left\llbracket
f^{(2)\flat}_{s}\left(
\mathfrak{Q}^{(2)}
\right)
\right\rrbracket_{\varphi(s)}.
\tag{5}
\end{align*}
\end{flushleft}

The case of $\mathfrak{P}^{(2)}$ being a second-order path of length strictly greater than one containing its first second-order echelon on its first step follows.

If $i>0$, that is, if $\mathfrak{P}^{(2)}$ is a second-order path of length strictly greater than one containing its first second-order echelon on a step different from the initial one, then according to Proposition~\ref{PDPthExt}, we have that
$$
f^{(2)\flat}_{s}\left(
\mathfrak{P}^{(2)}
\right)
=
f^{(2)\flat}_{s}\left(
\mathfrak{P}^{(2),i,\bb{\mathfrak{P}^{(2)}}-1}
\right)
\circ_{\varphi(s)}^{1\mathbf{Pth}_{\boldsymbol{\mathcal{B}}^{(2)}}}
f^{(2)\flat}_{s}\left(
\mathfrak{P}^{(2),0,i-1}
\right).
$$

Since $\mathrm{CH}^{(2)}_{\boldsymbol{\mathcal{A}}^{(2)}, s}(\mathfrak{P}^{(2)}) \in \eta^{(2,\mathcal{A}^{(2)})}[\mathcal{A}^{(2)}]_{s}^{\neg\mathrm{int}}$ and $(\mathfrak{P}^{(2)}, \mathfrak{Q}^{(2)}) \in \mathrm{Ker}(\mathrm{CH}^{(2)}_{\boldsymbol{\mathcal{A}}^{(2)}})_{s}$, we have, by Lemma~\ref{LDCHEchNInt}, that $\mathfrak{Q}^{(2)}$ is a second-order path of length strictly greater than one containing its first second-order echelon on a step different from the initial one.

Thus, according to Proposition~\ref{PDPthExt}, we have that
$$
f^{(2)\flat}_{s}\left(
\mathfrak{Q}^{(2)}
\right)
=
f^{(2)\flat}_{s}\left(
\mathfrak{Q}^{(2),i,\bb{\mathfrak{P}^{(2)}}-1}
\right)
\circ_{\varphi(s)}^{1\mathbf{Pth}_{\boldsymbol{\mathcal{B}}^{(2)}}}
f^{(2)\flat}_{s}\left(
\mathfrak{Q}^{(2),0,i-1}
\right).
$$

Since $(\mathfrak{P}^{(2)}, \mathfrak{Q}^{(2)}) \in \mathrm{Ker}(\mathrm{CH}^{(2)}_{\boldsymbol{\mathcal{A}}^{(2)}})_{s}$, we have that $(\mathfrak{P}^{(2),i,\bb{\mathfrak{P}^{(2)}}-1}, \mathfrak{Q}^{(2),i,\bb{\mathfrak{Q}^{(2)}}-1}$ and $(\mathfrak{P}^{(2),0,i-1}, \mathfrak{Q}^{(2),0,i-1})$ are in $\mathrm{Ker}(\mathrm{CH}^{(2)})_{s}$. Note that, according to Definition~\ref{DDOrd}, we have that $(\mathfrak{P}^{(2),0,i-1}, s)$ and $(\mathfrak{P}^{(2),i,\bb{\mathfrak{P}^{(2)}}-1}, s)$$\prec_{\mathbf{Pth}_{\boldsymbol{\mathcal{A}}^{(2)}}}$-precede $(\mathfrak{P}^{(2)}, s)$.

Therefore, by the inductive hypothesis, the second-order paths $\mathfrak{P}^{(2),0,i-1}$ and $\mathfrak{Q}^{(2),0,i-1}$, and the second-order paths $\mathfrak{P}^{(2),i,\bb{\mathfrak{P}^{(2)}}-1}$ and $\mathfrak{Q}^{(2),i,\bb{\mathfrak{Q}^{(2)}}-1}$ satisfy
\begin{align*}
\left\llbracket
f^{(2)\flat}_{s}\left(
\mathfrak{P}^{(2),0,i-1}
\right)
\right\rrbracket_{\varphi(s)}
&=
\left\llbracket
f^{(2)\flat}_{s}\left(
\mathfrak{Q}^{(2),0,i-1}
\right)
\right\rrbracket_{\varphi(s)}
\\
\left\llbracket
f^{(2)\flat}_{s}\left(
\mathfrak{P}^{(2),i,\bb{\mathfrak{P}^{(2)}}-1}
\right)
\right\rrbracket_{\varphi(s)}
&=
\left\llbracket
f^{(2)\flat}_{s}\left(
\mathfrak{Q}^{(2),i,\bb{\mathfrak{Q}^{(2)}}-1}
\right)
\right\rrbracket_{\varphi(s)}.
\end{align*}

Thus, the following chain of equalities holds
\begin{flushleft}
$
\left\llbracket
f^{(2)\flat}_{s}\left(
\mathfrak{P}^{(2)}
\right)
\right\rrbracket_{\varphi(s)}
$
\allowdisplaybreaks
\begin{align*}
&=
\left\llbracket
f^{(2)\flat}_{s}\left(
\mathfrak{P}^{(2),i,\bb{\mathfrak{P}^{(2)}}-1}
\right)
\circ_{\varphi(s)}^{1\mathbf{Pth}_{\boldsymbol{\mathcal{B}}^{(2)}}}
f^{(2)\flat}_{s}\left(
\mathfrak{P}^{(2),0,i-1}
\right)
\right\rrbracket_{\varphi(s)}
\tag{1}
\\
&=
\left\llbracket
f^{(2)\flat}_{s}\left(
\mathfrak{P}^{(2),i,\bb{\mathfrak{P}^{(2)}}-1}
\right)
\right\rrbracket_{\varphi(s)}
\circ_{\varphi(s)}^{1\llbracket\mathbf{Pth}_{\boldsymbol{\mathcal{B}}^{(2)}}\rrbracket}
\left\llbracket
f^{(2)\flat}_{s}\left(
\mathfrak{P}^{(2),0,i-1}
\right)
\right\rrbracket_{\varphi(s)}
\tag{2}
\\
&=
\left\llbracket
f^{(2)\flat}_{s}\left(
\mathfrak{Q}^{(2),i,\bb{\mathfrak{Q}^{(2)}}-1}
\right)
\right\rrbracket_{\varphi(s)}
\circ_{\varphi(s)}^{1\llbracket\mathbf{Pth}_{\boldsymbol{\mathcal{B}}^{(2)}}\rrbracket}
\left\llbracket
f^{(2)\flat}_{s}\left(
\mathfrak{Q}^{(2),0,i-1}
\right)
\right\rrbracket_{\varphi(s)}
\tag{3}
\\
&=
\left\llbracket
f^{(2)\flat}_{s}\left(
\mathfrak{Q}^{(2),i,\bb{\mathfrak{Q}^{(2)}}-1}
\right)
\circ_{\varphi(s)}^{1\mathbf{Pth}_{\boldsymbol{\mathcal{B}}^{(2)}}}
f^{(2)\flat}_{s}\left(
\mathfrak{Q}^{(2),0,i-1}
\right)
\right\rrbracket_{\varphi(s)}
\tag{4}
\\
&=
\left\llbracket
f^{(2)\flat}_{s}\left(
\mathfrak{Q}^{(2)}
\right)
\right\rrbracket_{\varphi(s)}
\tag{5}
\end{align*}
\end{flushleft}

The case of $\mathfrak{P}^{(2)}$ being a second-order path of length strictly greater than one containing its first second-order echelon on a step different from the initial one follows.

Case (1) follows.

If (2), i.e., if $\mathfrak{P}^{(2)}$ is an echelonless second-order path, it could be the case that (2.1) $\mathfrak{P}^{(2)}$ is an echelonless second-order path that is not head-constant, or (2.2) $\mathfrak{P}^{(2)}$ is a head-constant non-coherent echelonless second-order path or (2.3) $\mathfrak{P}^{(2)}$ is a head-constant coherent echelonless second-order path.

If $(2.1)$, let $i \in \bb{\mathfrak{P}^{(2)}}$ be the greatest index for which $\mathfrak{P}^{(2),0,i}$ is a head-constant second-order path. Then, according to Proposition~\ref{PDPthExt}, we have that
$$
f^{(2)\flat}_{s}\left(
\mathfrak{P}^{(2)}
\right)
=
f^{(2)\flat}_{s}\left(
\mathfrak{P}^{(2),i+1,\bb{\mathfrak{P}^{(2)}}-1}
\right)
\circ_{\varphi(s)}^{1\mathbf{Pth}_{\boldsymbol{\mathcal{B}}^{(2)}}}
f^{(2)\flat}_{s}\left(
\mathfrak{P}^{(2),0,i}
\right).
$$

Since $f^{(2)\flat}_{s}(\mathfrak{P}^{(2)}) \in \mathrm{T}_{\Sigma^{\boldsymbol{\mathcal{A}}^{(2)}}}(X)^{\neg\mathsf{HdC}}$ and $(\mathfrak{P}^{(2)}, \mathfrak{Q}^{(2)}) \in \mathrm{Ker}(\mathrm{CH}^{(2)}_{\boldsymbol{\mathcal{A}}^{(2)}})_{s}$, we have, by Lemma~\ref{LDCHNEchNHd}, that $\mathfrak{Q}^{(2)}$ is an echelonless second-order path that is not head-constant.

Thus, if $j \in \bb{\mathfrak{Q}^{(2)}}$ is the greatest index for which $\mathfrak{Q}^{(2),0,j}$ is a head-constant second-order path, according to Proposition~\ref{PDPthExt}, we have that
$$
f^{(2)\flat}_{s}\left(
\mathfrak{Q}^{(2)}
\right)
=
f^{(2)\flat}_{s}\left(
\mathfrak{Q}^{(2),j+1,\bb{\mathfrak{Q}^{(2)}}-1}
\right)
\circ_{\varphi(s)}^{1\mathbf{Pth}_{\boldsymbol{\mathcal{B}}^{(2)}}}
f^{(2)\flat}_{s}\left(
\mathfrak{Q}^{(2),0,j}
\right).
$$

Since $(\mathfrak{P}^{(2)}, \mathfrak{Q}^{(2)}) \in \mathrm{Ker}(\mathrm{CH}^{(2)}_{\boldsymbol{\mathcal{A}}^{(2)}})_{s}$, we have that $(\mathfrak{P}^{(2),i+1,\bb{\mathfrak{P}^{(2)}}-1},\mathfrak{Q}^{(2),j+1,\bb{\mathfrak{Q}^{(2)}}-1})$ and $(\mathfrak{P}^{(2),0,i},\mathfrak{Q}^{(2),0,j})$ are in $\mathrm{Ker}(\mathrm{CH}^{(2)}_{\boldsymbol{\mathcal{A}}^{(2)}})_{s}$. Note that, according to Definition~\ref{DDOrd}, we have that $(\mathfrak{P}^{(2),0,i}, s)$ and $(\mathfrak{P}^{(2),i+1,\bb{\mathfrak{P}^{(2)}}-1})$$\prec_{\mathbf{Pth}_{\boldsymbol{\mathcal{A}}^{(2)}}}$-precede $(\mathfrak{P}^{(2)},s)$.

Therefore, by the inductive hypothesis, the second-order paths $\mathfrak{P}^{(2),0,i}$ and $\mathfrak{Q}^{(2),0,j}$, and the second-order paths $\mathfrak{P}^{(2),i+1,\bb{\mathfrak{P}^{(2)}}-1}$ and $\mathfrak{Q}^{(2),j+1,\bb{\mathfrak{Q}^{(2)}}-1}$ satisfy
\begin{align*}
\left\llbracket
f^{(2)\flat}_{s}\left(
\mathfrak{P}^{(2),0,i}
\right)
\right\rrbracket_{\varphi(s)}
&=
\left\llbracket
f^{(2)\flat}_{s}\left(
\mathfrak{Q}^{(2),0,j}
\right)
\right\rrbracket_{\varphi(s)}
\\
\left\llbracket
f^{(2)\flat}_{s}\left(
\mathfrak{P}^{(2),i+1,\bb{\mathfrak{P}^{(2)}}-1}
\right)
\right\rrbracket_{\varphi(s)}
&=
\left\llbracket
f^{(2)\flat}_{s}\left(
\mathfrak{Q}^{(2),j+1,\bb{\mathfrak{Q}^{(2)}}-1}
\right)
\right\rrbracket_{\varphi(s)}.
\end{align*}

Thus, the following chain of equalities holds
\begin{flushleft}
$
\left\llbracket
f^{(2)\flat}_{s}\left(
\mathfrak{P}^{(2)}
\right)
\right\rrbracket_{\varphi(s)}
$
\allowdisplaybreaks
\begin{align*}
&=
\left\llbracket
f^{(2)\flat}_{s}\left(
\mathfrak{P}^{(2),i+1,\bb{\mathfrak{P}^{(2)}}-1}
\right)
\circ_{\varphi(s)}^{1\mathbf{Pth}_{\boldsymbol{\mathcal{B}}^{(2)}}}
f^{(2)\flat}_{s}\left(
\mathfrak{P}^{(2),0,i}
\right)
\right\rrbracket_{\varphi(s)}
\tag{1}
\\
&=
\left\llbracket
f^{(2)\flat}_{s}\left(
\mathfrak{P}^{(2),i+1,\bb{\mathfrak{P}^{(2)}}-1}
\right)
\right\rrbracket_{\varphi(s)}
\circ_{\varphi(s)}^{1\llbracket\mathbf{Pth}_{\boldsymbol{\mathcal{B}}^{(2)}}\rrbracket}
\left\llbracket
f^{(2)\flat}_{s}\left(
\mathfrak{P}^{(2),0,i}
\right)
\right\rrbracket_{\varphi(s)}
\tag{2}
\\
&=
\left\llbracket
f^{(2)\flat}_{s}\left(
\mathfrak{Q}^{(2),j+1,\bb{\mathfrak{Q}^{(2)}}-1}
\right)
\right\rrbracket_{\varphi(s)}
\circ_{\varphi(s)}^{1\llbracket\mathbf{Pth}_{\boldsymbol{\mathcal{B}}^{(2)}}\rrbracket}
\left\llbracket
f^{(2)\flat}_{s}\left(
\mathfrak{Q}^{(2),0,j}
\right)
\right\rrbracket_{\varphi(s)}
\tag{3}
\\
&=
\left\llbracket
f^{(2)\flat}_{s}\left(
\mathfrak{Q}^{(2),j+1,\bb{\mathfrak{Q}^{(2)}}-1}
\right)
\circ_{\varphi(s)}^{1\mathbf{Pth}_{\boldsymbol{\mathcal{B}}^{(2)}}}
f^{(2)\flat}_{s}\left(
\mathfrak{Q}^{(2),0,j}
\right)
\right\rrbracket_{\varphi(s)}
\tag{4}
\\
&=
\left\llbracket
f^{(2)\flat}_{s}\left(
\mathfrak{Q}^{(2)}
\right)
\right\rrbracket_{\varphi(s)}
\tag{5}
\end{align*}
\end{flushleft}

The case of $\mathfrak{P}^{(2)}$ being an echelonless second-order path that is not head-
constant follows.

If $(2.2)$, let $i \in \bb{\mathfrak{P}^{(2)}}$ be the greatest index for which $\mathfrak{P}^{(2),0,i}$ is a coherent second-order path. Then, according to Proposition~\ref{PDPthExt}, we have that
$$
f^{(2)\flat}_{s}\left(
\mathfrak{P}^{(2)}
\right)
=
f^{(2)\flat}_{s}\left(
\mathfrak{P}^{(2),i+1,\bb{\mathfrak{P}^{(2)}}-1}
\right)
\circ_{\varphi(s)}^{1\mathbf{Pth}_{\boldsymbol{\mathcal{B}}^{(2)}}}
f^{(2)\flat}_{s}\left(
\mathfrak{P}^{(2),0,i}
\right).
$$

Since $f^{(2)\flat}_{s}(\mathfrak{P}^{(2)}) \in \mathrm{T}_{\Sigma^{\boldsymbol{\mathcal{A}}^{(2)}}}(X)^{\mathsf{HdC}\&\neg\mathsf{C}}$ and $(\mathfrak{P}^{(2)}, \mathfrak{Q}^{(2)}) \in \mathrm{Ker}(\mathrm{CH}^{(2)}_{\boldsymbol{\mathcal{A}}^{(2)}})_{s}$, we have, by Lemma~\ref{LDCHNEchHdNC}, that $\mathfrak{Q}^{(2)}$ is an echelonless second-order path that is a head-constant non-coherent echelonless second-order path.

Thus, if $j \in \bb{\mathfrak{Q}^{(2)}}$ is the greatest index for which $\mathfrak{Q}^{(2),0,j}$ is a coherent second-order path, according to Proposition~\ref{PDPthExt}, we have that
$$
f^{(2)\flat}_{s}\left(
\mathfrak{Q}^{(2)}
\right)
=
f^{(2)\flat}_{s}\left(
\mathfrak{Q}^{(2),j+1,\bb{\mathfrak{Q}^{(2)}}-1}
\right)
\circ_{\varphi(s)}^{1\mathbf{Pth}_{\boldsymbol{\mathcal{B}}^{(2)}}}
f^{(2)\flat}_{s}\left(
\mathfrak{Q}^{(2),0,j}
\right).
$$

Since $(\mathfrak{P}^{(2)}, \mathfrak{Q}^{(2)}) \in \mathrm{Ker}(\mathrm{CH}^{(2)}_{\boldsymbol{\mathcal{A}}^{(2)}})_{s}$, we have that $(\mathfrak{P}^{(2),i+1,\bb{\mathfrak{P}^{(2)}}-1},\mathfrak{Q}^{(2),j+1,\bb{\mathfrak{Q}^{(2)}}-1})$ and $(\mathfrak{P}^{(2),0,i},\mathfrak{Q}^{(2),0,j})$ are in $\mathrm{Ker}(\mathrm{CH}^{(2)}_{\boldsymbol{\mathcal{A}}^{(2)}})_{s}$. Note that, according to Definition~\ref{DDOrd}, we have that $(\mathfrak{P}^{(2),0,i}, s)$ and $(\mathfrak{P}^{(2),i+1,\bb{\mathfrak{P}^{(2)}}-1})$$\prec_{\mathbf{Pth}_{\boldsymbol{\mathcal{A}}^{(2)}}}$-precede $(\mathfrak{P}^{(2)},s)$.

Therefore, by the inductive hypothesis, the second-order paths $\mathfrak{P}^{(2),0,i}$ and $\mathfrak{Q}^{(2),0,j}$, and the second-order paths $\mathfrak{P}^{(2),i+1,\bb{\mathfrak{P}^{(2)}}-1}$ and $\mathfrak{Q}^{(2),j+1,\bb{\mathfrak{Q}^{(2)}}-1}$ satisfy
\begin{align*}
\left\llbracket
f^{(2)\flat}_{s}\left(
\mathfrak{P}^{(2),0,i}
\right)
\right\rrbracket_{\varphi(s)}
&=
\left\llbracket
f^{(2)\flat}_{s}\left(
\mathfrak{Q}^{(2),0,j}
\right)
\right\rrbracket_{\varphi(s)}
\\
\left\llbracket
f^{(2)\flat}_{s}\left(
\mathfrak{P}^{(2),i+1,\bb{\mathfrak{P}^{(2)}}-1}
\right)
\right\rrbracket_{\varphi(s)}
&=
\left\llbracket
f^{(2)\flat}_{s}\left(
\mathfrak{Q}^{(2),j+1,\bb{\mathfrak{Q}^{(2)}}-1}
\right)
\right\rrbracket_{\varphi(s)}.
\end{align*}

Thus, the following chain of equalities holds
\begin{flushleft}
$
\left\llbracket
f^{(2)\flat}_{s}\left(
\mathfrak{P}^{(2)}
\right)
\right\rrbracket_{\varphi(s)}
$
\allowdisplaybreaks
\begin{align*}
&=
\left\llbracket
f^{(2)\flat}_{s}\left(
\mathfrak{P}^{(2),i+1,\bb{\mathfrak{P}^{(2)}}-1}
\right)
\circ_{\varphi(s)}^{1\mathbf{Pth}_{\boldsymbol{\mathcal{B}}^{(2)}}}
f^{(2)\flat}_{s}\left(
\mathfrak{P}^{(2),0,i}
\right)
\right\rrbracket_{\varphi(s)}
\tag{1}
\\
&=
\left\llbracket
f^{(2)\flat}_{s}\left(
\mathfrak{P}^{(2),i+1,\bb{\mathfrak{P}^{(2)}}-1}
\right)
\right\rrbracket_{\varphi(s)}
\circ_{\varphi(s)}^{1\llbracket\mathbf{Pth}_{\boldsymbol{\mathcal{B}}^{(2)}}\rrbracket}
\left\llbracket
f^{(2)\flat}_{s}\left(
\mathfrak{P}^{(2),0,i}
\right)
\right\rrbracket_{\varphi(s)}
\tag{2}
\\
&=
\left\llbracket
f^{(2)\flat}_{s}\left(
\mathfrak{Q}^{(2),j+1,\bb{\mathfrak{Q}^{(2)}}-1}
\right)
\right\rrbracket_{\varphi(s)}
\circ_{\varphi(s)}^{1\llbracket\mathbf{Pth}_{\boldsymbol{\mathcal{B}}^{(2)}}\rrbracket}
\left\llbracket
f^{(2)\flat}_{s}\left(
\mathfrak{Q}^{(2),0,j}
\right)
\right\rrbracket_{\varphi(s)}
\tag{3}
\\
&=
\left\llbracket
f^{(2)\flat}_{s}\left(
\mathfrak{Q}^{(2),j+1,\bb{\mathfrak{Q}^{(2)}}-1}
\right)
\circ_{\varphi(s)}^{1\mathbf{Pth}_{\boldsymbol{\mathcal{B}}^{(2)}}}
f^{(2)\flat}_{s}\left(
\mathfrak{Q}^{(2),0,j}
\right)
\right\rrbracket_{\varphi(s)}
\tag{4}
\\
&=
\left\llbracket
f^{(2)\flat}_{s}\left(
\mathfrak{Q}^{(2)}
\right)
\right\rrbracket_{\varphi(s)}
\tag{5}
\end{align*}
\end{flushleft}

The case of $\mathfrak{P}^{(2)}$ being a head-constant echelonless second-order path that is
not coherent follows.

If $(2.3)$, then there exists a unique word $\mathbf{s} \in S^{\star} - \{\lambda\}$ and a unique operation symbol $\tau \in \Sigma^{\boldsymbol{\mathcal{A}}^{(2)}}_{\mathbf{s}, s}$ associated to $\mathfrak{P}^{(2)}$. Let $(\mathfrak{P}^{(2)}_{j})_{j \in \bb{\mathbf{s}}}$ be the family of second-order paths in $\mathbf{Pth}_{\boldsymbol{\mathcal{A}}^{(2)}, \mathbf{s}}$ which, in virtue of Lemma~\ref{LDPthExtract}, we can extract from $\mathfrak{P}^{(2)}$. Then, according to Proposition~\ref{PDPthExt}, we have that
$$
f^{(2)\flat}_{s}\left(
\mathfrak{P}^{(2)}
\right)
=
\tau^{\mathbf{Pth}_{\boldsymbol{\mathcal{B}}^{(2)}}^{\mathbf{f}^{(2)}}}\left(\left(
f^{(2)\flat}_{s_{j}}\left(
\mathfrak{P}^{(2)}_{j}
\right)
\right)_{j\in\bb{\mathbf{s}}}\right).
$$

Since $\mathrm{CH}^{(2)}_{\boldsymbol{\mathcal{A}}^{(2)}, s}(\mathfrak{P}^{(2)}) \in \mathcal{T}(\tau, \mathrm{T}_{\Sigma^{\boldsymbol{\mathcal{A}}^{(2)}}(X)})_{1}$, which is a subset of $\mathrm{T}_{\Sigma^{\boldsymbol{\mathcal{A}}^{(2)}}}(X)_{s}^{\mathsf{HdC}\&\mathsf{C}}$, and $(\mathfrak{P}^{(2)}, \mathfrak{Q}^{(2)}) \in \mathrm{Ker}(\mathrm{CH}^{(2)}_{\boldsymbol{\mathcal{A}}^{(2)}})_{s}$ we have, by Lemma~\ref{LDCHNEchHdC}, that $\mathfrak{Q}^{(2)}$ is a head-constant coherent echelonless second-order path associated to $\tau$, the same operation symbol as that associated to $\mathfrak{P}^{(2)}$.

Let $(\mathfrak{Q}^{(2)}_{j})_{j\in\bb{\mathbf{s}}}$ be the family of second-order paths in $\mathrm{Pth}_{\boldsymbol{\mathcal{A}}^{(2)}, \mathbf{s}}$ which, by Lemma~\ref{LDPthExtract}, we can extract from $\mathfrak{Q}^{(2)}$. Then, according to Proposition~\ref{PDPthExt}, we have that 
$$
f^{(2)\flat}_{s}\left(
\mathfrak{Q}^{(2)}
\right)
=
\tau^{\mathbf{Pth}_{\boldsymbol{\mathcal{B}}^{(2)}}^{\mathbf{f}^{(2)}}}\left(\left(
f^{(2)\flat}_{s_{j}}\left(
\mathfrak{Q}^{(2)}_{j}
\right)
\right)_{j\in\bb{\mathbf{s}}}\right).
$$

Since $(\mathfrak{P}^{(2)}, \mathfrak{Q}^{(2)}) \in \mathrm{Ker}(\mathrm{CH}^{(2)}_{\boldsymbol{\mathcal{A}}^{(2)}})_{s}$, we have, for every $j \in \bb{\mathbf{s}}$, that $(\mathfrak{P}^{(2)}_{j}, \mathfrak{Q}^{(2)}_{j}) \in \mathrm{Ker}(\mathrm{CH}^{(2)}_{\boldsymbol{\mathcal{A}}^{(2)}})_{s_{j}}$. Note that, accoording to Definition~\ref{DDOrd}, we have that, for every $j \in \bb{\mathbf{s}}$, $(\mathfrak{P}^{(2)}_{j}, s_{j})$$\prec_{\mathbf{Pth}_{\boldsymbol{\mathcal{A}}^{(2)}}}$-precedes $(\mathfrak{P}^{(2)},s)$.

Therefore, by the inductive hypothesis, for every $j \in \bb{\mathbf{s}}$, the second-order paths $\mathfrak{P}^{(2)}_{j}$ and $\mathfrak{Q}^{(2)}_{j}$ satisfy
$$
\left\llbracket
f^{(2)\flat}_{s_{j}}\left(
\mathfrak{P}^{(2)}_{j}
\right)
\right\rrbracket_{\varphi(s_{j})}
=
\left\llbracket
f^{(2)\flat}_{s_{j}}\left(
\mathfrak{Q}^{(2)}_{j}
\right)
\right\rrbracket_{\varphi(s_{j})}
$$

Thus, the following chain of equalities holds
\allowdisplaybreaks
\begin{align*}
\left\llbracket
f^{(2)\flat}_{s}\left(
\mathfrak{P}^{(2)}
\right)
\right\rrbracket_{\varphi(s)}
&=
\left\llbracket
\tau^{\mathbf{Pth}_{\boldsymbol{\mathcal{B}}^{(2)}}^{\mathbf{f}^{(2)}}}\left(\left(
f^{(2)\flat}_{s_{j}}\left(
\mathfrak{P}^{(2)}_{j}
\right)
\right)_{j\in\bb{\mathbf{s}}}\right)
\right\rrbracket_{\varphi(s)}
\tag{1}
\\
&=
\tau^{\llbracket\mathbf{Pth}_{\boldsymbol{\mathcal{B}}^{(2)}}^{\mathbf{f}^{(2)}}\rrbracket}\left(\left(
\left\llbracket
f^{(2)\flat}_{s_{j}}\left(
\mathfrak{P}^{(2)}_{j}
\right)
\right\rrbracket_{\varphi(s_{j})}
\right)_{j\in\bb{\mathbf{s}}}\right)
\tag{2}
\\
&=
\tau^{\llbracket\mathbf{Pth}_{\boldsymbol{\mathcal{B}}^{(2)}}^{\mathbf{f}^{(2)}}\rrbracket}\left(\left(
\left\llbracket
f^{(2)\flat}_{s_{j}}\left(
\mathfrak{Q}^{(2)}_{j}
\right)
\right\rrbracket_{\varphi(s_{j})}
\right)_{j\in\bb{\mathbf{s}}}\right)
\tag{3}
\\
&=
\left\llbracket
\tau^{\mathbf{Pth}_{\boldsymbol{\mathcal{B}}^{(2)}}^{\mathbf{f}^{(2)}}}\left(\left(
f^{(2)\flat}_{s_{j}}\left(
\mathfrak{Q}^{(2)}_{j}
\right)
\right)_{j\in\bb{\mathbf{s}}}\right)
\right\rrbracket_{\varphi(s)}
\tag{4}
\\
&=
\left\llbracket
f^{(2)\flat}_{s}\left(
\mathfrak{Q}^{(2)}
\right)
\right\rrbracket_{\varphi(s)}
\tag{5}
\end{align*}

The case of $\mathfrak{P}^{(2)}$ being a head-constant coherent echelonless second-order path
follows.

Case (1) follows.

This proves that $\mathrm{Ker}(\mathrm{CH}^{(2)}_{\boldsymbol{\mathcal{A}}^{(2)}})$ is included in $\mathrm{Ker}(\mathrm{pr}_{\boldsymbol{\mathcal{B}}^{(2)}, \varphi}^{\llbracket\cdot\rrbracket} \circ f^{(2)\flat})$.

{\sffamily (ii) $\Upsilon^{[1]} \subseteq \mathrm{Ker}(\mathrm{pr}_{\boldsymbol{\mathcal{B}}^{(2)}, \varphi}^{\llbracket\cdot\rrbracket} \circ f^{(2)\flat})$.}

Let us recall from Definition~\ref{DDUpsCong} that $\Upsilon^{[1]}$ is defined as the samllest $\Sigma^{\boldsymbol{\mathcal{A}}^{(2)}}$-congruence containing the relation $\Upsilon^{(1)}$, introduced in Definition~\ref{DDUps}. Therefore, since we have already proven that $\mathrm{pr}_{\boldsymbol{\mathcal{B}}^{(2)}, \varphi}^{\llbracket\cdot\rrbracket} \circ f^{(2)\flat}$ is a $\Sigma^{\boldsymbol{\mathcal{A}}^{(2)}}$-homomorphism, in order to check that $\Upsilon^{[1]}$ is included in $\mathrm{Ker}(\mathrm{pr}_{\boldsymbol{\mathcal{B}}^{(2)}, \varphi}^{\llbracket\cdot\rrbracket} \circ f^{(2)\flat})$ it suffices to prove that $\Upsilon^{(1)}$ is included in $\mathrm{Ker}(\mathrm{pr}_{\boldsymbol{\mathcal{B}}^{(2)}, \varphi}^{\llbracket\cdot\rrbracket} \circ f^{(2)\flat})$.

Following Definition~\ref{DDUps}, we will consider the different cases defining $\Upsilon^{(1)}$ individually.

(1)
For every sort $s \in S$ and every second-order path $\mathfrak{P}^{(2)}$ in $\mathrm{Pth}_{\boldsymbol{\mathcal{A}}^{(2)}, s}$,
$$
\left(
\mathfrak{P}^{(2)},
\mathfrak{P}^{(2)} \circ_{s}^{0\mathbf{Pth}_{\boldsymbol{\mathcal{A}}^{(2)}}} \mathrm{sc}_{s}^{0\mathbf{Pth}_{\boldsymbol{\mathcal{A}}^{(2)}}}\left(\mathfrak{P}^{(2)}\right)
\right)
\in
\Upsilon_{s}^{(1)}.
$$

We want to check that
$$
\left\llbracket
f^{(2)\flat}_{s}\left(
\mathfrak{P}^{(2)}
\right)
\right\rrbracket_{\varphi(s)}
\\
=
\left\llbracket
f^{(2)\flat}_{s}\left(
\mathfrak{P}^{(2)} \circ_{s}^{0\mathbf{Pth}_{\boldsymbol{\mathcal{A}}^{(2)}}} \mathrm{sc}_{s}^{0\mathbf{Pth}_{\boldsymbol{\mathcal{A}}^{(2)}}}\left(\mathfrak{P}^{(2)}\right)
\right)
\right\rrbracket_{\varphi(s)}.
$$

The following chain of equalities holds
\begin{flushleft}
$
\left\llbracket
f^{(2)\flat}_{s}\left(
\mathfrak{P}^{(2)} \circ_{s}^{0\mathbf{Pth}_{\boldsymbol{\mathcal{A}}^{(2)}}} \mathrm{sc}_{s}^{0\mathbf{Pth}_{\boldsymbol{\mathcal{A}}^{(2)}}}\left(\mathfrak{P}^{(2)}\right)
\right)
\right\rrbracket_{\varphi(s)}
$
\allowdisplaybreaks
\begin{align*}
&=
\left\llbracket
f^{(2)\flat}_{s}\left(
\mathfrak{P}^{(2)}
\right)
\right\rrbracket_{\varphi(s)}
\circ_{s}^{0\llbracket\mathbf{Pth}_{\boldsymbol{\mathcal{B}}^{(2)}}^{\mathbf{f}^{(2)}}\rrbracket} 
\mathrm{sc}_{s}^{0\llbracket\mathbf{Pth}_{\boldsymbol{\mathcal{B}}^{(2)}}^{\mathbf{f}^{(2)}}\rrbracket}\left(
\left\llbracket
f^{(2)\flat}_{s}\left(
\mathfrak{P}^{(2)}
\right)
\right\rrbracket_{\varphi(s)}
\right)
\tag{1}
\\
&=
\left\llbracket
f^{(2)\flat}_{s}\left(
\mathfrak{P}^{(2)}
\right)
\right\rrbracket_{\varphi(s)}
\circ_{\varphi(s)}^{0\llbracket\mathbf{Pth}_{\boldsymbol{\mathcal{B}}^{(2)}}\rrbracket} 
\mathrm{sc}_{\varphi(s)}^{0\llbracket\mathbf{Pth}_{\boldsymbol{\mathcal{B}}^{(2)}}\rrbracket}\left(
\left\llbracket
f^{(2)\flat}_{s}\left(
\mathfrak{P}^{(2)}
\right)
\right\rrbracket_{\varphi(s)}
\right)
\tag{2}
\\
&=
\left\llbracket
f^{(2)\flat}_{s}\left(
\mathfrak{P}^{(2)}
\right)
\right\rrbracket_{\varphi(s)}.
\tag{3}
\end{align*}
\end{flushleft}

The first equality follows from the fact that $\mathrm{pr}_{\boldsymbol{\mathcal{B}}^{(2)}, \varphi}^{\llbracket\cdot\rrbracket} \circ f^{(2)\flat}$ is a $\Sigma^{\boldsymbol{\mathcal{A}}^{(2)}}$-homomorphism;
the second equality unravels the definition of the operation symbols $\mathrm{sc}^{0}_{s}$ and $\circ^{0}_{s}$ in the partial $\Sigma^{\boldsymbol{\mathcal{A}}^{(2)}}$-algebra $\llbracket \mathbf{Pth}_{\boldsymbol{\mathcal{B}}^{(2)}}^{\mathbf{f}^{(2)}} \rrbracket$ introduced in Proposition~\ref{PDQPthBCatAlg};
finally, the last equality follows from Axiom~(A5) since, by Proposition~\ref{PDPthVar}, $\llbracket\mathbf{Pth}_{\boldsymbol{\mathcal{B}}^{(2)}}\rrbracket$ belongs to $\mathbf{PAlg}(\boldsymbol{\mathcal{E}}^{\boldsymbol{\mathcal{B}}^{(2)}})$.

This proves Case (1).

(2) 
For every sort $s \in S$ and every second-order path $\mathfrak{P}^{(2)}$ in $\mathrm{Pth}_{\boldsymbol{\mathcal{A}}^{(2)}, s}$,
$$
\left(
\mathfrak{P}^{(2)},
\mathrm{tg}_{s}^{0\mathbf{Pth}_{\boldsymbol{\mathcal{A}}^{(2)}}}\left(\mathfrak{P}^{(2)}\right) \circ_{s}^{0\mathbf{Pth}_{\boldsymbol{\mathcal{A}}^{(2)}}} \mathfrak{P}^{(2)}
\right)
\in
\Upsilon_{s}^{(1)}.
$$

We want to check that
$$
\left\llbracket
f^{(2)\flat}_{s}\left(
\mathfrak{P}^{(2)}
\right)
\right\rrbracket_{\varphi(s)}
\\
=
\left\llbracket
f^{(2)\flat}_{s}\left(
\mathrm{tg}_{s}^{0\mathbf{Pth}_{\boldsymbol{\mathcal{A}}^{(2)}}}\left(\mathfrak{P}^{(2)}\right) \circ_{s}^{0\mathbf{Pth}_{\boldsymbol{\mathcal{A}}^{(2)}}} \mathfrak{P}^{(2)}
\right)
\right\rrbracket_{\varphi(s)}.
$$

This case follows by a similar argument to that presented for the Case (1).

This proves Case (2).

(3)
For every sort $s \in S$ and every second-order paths $\mathfrak{P}^{(2)}$, $\mathfrak{Q}^{(2)}$ and $\mathfrak{R}^{(2)}$ in $\mathrm{Pth}_{\boldsymbol{\mathcal{A}}^{(2)}, s}$, satisfying that
\begin{align*}
\mathrm{sc}^{(0,2)}_{\boldsymbol{\mathcal{A}}^{(2)}, s}\left(
\mathfrak{R}^{(2)}
\right)
&=
\mathrm{tg}^{(0,2)}_{\boldsymbol{\mathcal{A}}^{(2)}, s}\left(
\mathfrak{Q}^{(2)}
\right);&
\mathrm{sc}^{(0,2)}_{\boldsymbol{\mathcal{A}}^{(2)}, s}\left(
\mathfrak{Q}^{(2)}
\right)
&=
\mathrm{tg}^{(0,2)}_{\boldsymbol{\mathcal{A}}^{(2)}, s}\left(
\mathfrak{P}^{(2)}
\right);
\end{align*}
then 
\begin{multline*}
\left(
\mathfrak{R}^{(2)} \circ^{0\mathbf{Pth}_{\boldsymbol{\mathcal{A}}^{(2)}}}_{s} \left( \mathfrak{Q}^{(2)} \circ^{0\mathbf{Pth}_{\boldsymbol{\mathcal{A}}^{(2)}}}_{s} \mathfrak{P}^{(2)} \right),
\right.
\\
\left.
\left( \mathfrak{R}^{(2)} \circ^{0\mathbf{Pth}_{\boldsymbol{\mathcal{A}}^{(2)}}}_{s} \mathfrak{Q}^{(2)} \right) \circ^{0\mathbf{Pth}_{\boldsymbol{\mathcal{A}}^{(2)}}}_{s} \mathfrak{P}^{(2)}
\right) \in \Upsilon_{s}^{(1)}.
\end{multline*}

We want to check that
\begin{multline*}
\left\llbracket
f^{(2)\flat}_{s} \left(
\mathfrak{R}^{(2)} \circ^{0\mathbf{Pth}_{\boldsymbol{\mathcal{A}}^{(2)}}}_{s} \left( \mathfrak{Q}^{(2)} \circ^{0\mathbf{Pth}_{\boldsymbol{\mathcal{A}}^{(2)}}}_{s} \mathfrak{P}^{(2)} \right)
\right)
\right\rrbracket_{\varphi(s)}
\\
=
\left\llbracket
f^{(2)\flat}_{s} \left(
\left( \mathfrak{R}^{(2)} \circ^{0\mathbf{Pth}_{\boldsymbol{\mathcal{A}}^{(2)}}}_{s} \mathfrak{Q}^{(2)} \right) \circ^{0\mathbf{Pth}_{\boldsymbol{\mathcal{A}}^{(2)}}}_{s} \mathfrak{P}^{(2)}
\right)
\right\rrbracket_{\varphi(s)}.
\end{multline*}

The following chain of equalities holds
\begin{flushleft}
$
\left\llbracket
f^{(2)\flat}_{s} \left(
\mathfrak{R}^{(2)} \circ^{0\mathbf{Pth}_{\boldsymbol{\mathcal{A}}^{(2)}}}_{s} \left( \mathfrak{Q}^{(2)} \circ^{0\mathbf{Pth}_{\boldsymbol{\mathcal{A}}^{(2)}}}_{s} \mathfrak{P}^{(2)} \right)
\right)
\right\rrbracket_{\varphi(s)}
$
\allowdisplaybreaks
\begin{align*}
&=
\left\llbracket
f^{(2)\flat}_{s} \left(
\mathfrak{R}^{(2)}
\right)
\right\rrbracket_{\varphi(s)}
\circ^{0\llbracket\mathbf{Pth}_{\boldsymbol{\mathcal{B}}^{(2)}}^{\mathbf{f}^{(2)}}\rrbracket}_{s}
\\
&\hspace{3cm}
\left(
\left\llbracket
f^{(2)\flat}_{s} \left(
\mathfrak{Q}^{(2)} 
\right)
\right\rrbracket_{\varphi(s)}
\circ^{0\llbracket\mathbf{Pth}_{\boldsymbol{\mathcal{B}}^{(2)}}^{\mathbf{f}^{(2)}}\rrbracket}_{\varphi(s)}
\left\llbracket
f^{(2)\flat}_{s} \left(
\mathfrak{P}^{(2)}
\right)
\right\rrbracket_{\varphi(s)}
\right)
\tag{1}
\\
&=
\left\llbracket
f^{(2)\flat}_{s} \left(
\mathfrak{R}^{(2)}
\right)
\right\rrbracket_{\varphi(s)}
\circ^{0\llbracket\mathbf{Pth}_{\boldsymbol{\mathcal{B}}^{(2)}}\rrbracket}_{\varphi(s)}
\\
&\hspace{3cm}
\left(
\left\llbracket
f^{(2)\flat}_{s} \left(
\mathfrak{Q}^{(2)} 
\right)
\right\rrbracket_{\varphi(s)}
\circ^{0\llbracket\mathbf{Pth}_{\boldsymbol{\mathcal{B}}^{(2)}}\rrbracket}_{\varphi(s)}
\left\llbracket
f^{(2)\flat}_{s} \left(
\mathfrak{P}^{(2)}
\right)
\right\rrbracket_{\varphi(s)}
\right)
\tag{2}
\\
&=
\left(
\left\llbracket
f^{(2)\flat}_{s} \left(
\mathfrak{R}^{(2)}
\right)
\right\rrbracket_{\varphi(s)}
\circ^{0\llbracket\mathbf{Pth}_{\boldsymbol{\mathcal{B}}^{(2)}}\rrbracket}_{\varphi(s)}
\left\llbracket
f^{(2)\flat}_{s} \left(
\mathfrak{Q}^{(2)} 
\right)
\right\rrbracket_{\varphi(s)}
\right)
\\
&\hspace{7cm}
\circ^{0\llbracket\mathbf{Pth}_{\boldsymbol{\mathcal{B}}^{(2)}}\rrbracket}_{\varphi(s)}
\left\llbracket
f^{(2)\flat}_{s} \left(
\mathfrak{P}^{(2)}
\right)
\right\rrbracket_{\varphi(s)}
\tag{3}
\\
&=
\left(
\left\llbracket
f^{(2)\flat}_{s} \left(
\mathfrak{R}^{(2)}
\right)
\right\rrbracket_{\varphi(s)}
\circ^{0\llbracket\mathbf{Pth}_{\boldsymbol{\mathcal{B}}^{(2)}}^{\mathbf{f}^{(2)}}\rrbracket}_{s}
\left\llbracket
f^{(2)\flat}_{s} \left(
\mathfrak{Q}^{(2)} 
\right)
\right\rrbracket_{\varphi(s)}
\right)
\\
&\hspace{7cm}
\circ^{0\llbracket\mathbf{Pth}_{\boldsymbol{\mathcal{B}}^{(2)}}^{\mathbf{f}^{(2)}}\rrbracket}_{s}
\left\llbracket
f^{(2)\flat}_{s} \left(
\mathfrak{P}^{(2)}
\right)
\right\rrbracket_{\varphi(s)}
\tag{4}
\\
&=
\left\llbracket
f^{(2)\flat}_{s} \left(
\left( \mathfrak{R}^{(2)} \circ^{0\mathbf{Pth}_{\boldsymbol{\mathcal{A}}^{(2)}}}_{s} \mathfrak{Q}^{(2)} \right) \circ^{0\mathbf{Pth}_{\boldsymbol{\mathcal{A}}^{(2)}}}_{s} \mathfrak{P}^{(2)}
\right)
\right\rrbracket_{\varphi(s)}.
\tag{5}
\end{align*}
\end{flushleft}

The first equality follows from the fact that $\mathrm{pr}_{\boldsymbol{\mathcal{B}}^{(2)}, \varphi}^{\llbracket\cdot\rrbracket} \circ f^{(2)\flat}$ is a $\Sigma^{\boldsymbol{\mathcal{A}}^{(2)}}$-homomorphism;
the second equality unravels the definition of the operation symbol $\circ^{0}_{s}$ in the partial $\Sigma^{\boldsymbol{\mathcal{A}}^{(2)}}$-algebra $\llbracket \mathbf{Pth}_{\boldsymbol{\mathcal{B}}^{(2)}}^{\mathbf{f}^{(2)}} \rrbracket$ introduced in Proposition~\ref{PDQPthBCatAlg};
the third equality follows from Axiom~(A6) since, by Proposition~\ref{PDPthVar}, $\llbracket\mathbf{Pth}_{\boldsymbol{\mathcal{B}}^{(2)}}\rrbracket$ belongs to $\mathbf{PAlg}(\boldsymbol{\mathcal{E}}^{\boldsymbol{\mathcal{B}}^{(2)}})$.
the fourth equality recovers the definition of the operation symbol $\circ^{0}_{s}$ in the partial $\Sigma^{\boldsymbol{\mathcal{A}}^{(2)}}$-algebra $\llbracket \mathbf{Pth}_{\boldsymbol{\mathcal{B}}^{(2)}}^{\mathbf{f}^{(2)}} \rrbracket$ introduced in Proposition~\ref{PDQPthBCatAlg};
finally, the last equality follows from the fact that $\mathrm{pr}_{\boldsymbol{\mathcal{B}}^{(2)}, \varphi}^{\llbracket\cdot\rrbracket} \circ f^{(2)\flat}$ is a $\Sigma^{\boldsymbol{\mathcal{A}}^{(2)}}$-homomorphism.

This proves Case (3).

(4)
For every word $\mathbf{s}$ and $s \in S$, for every operation symbol $\sigma \in \Sigma_{\mathbf{s}, s}$, for every pair of families of second-order paths $(\mathfrak{P}^{(2)}_{j})_{j \in \bb{\mathbf{s}}}$ and $(\mathfrak{Q}^{(2)}_{j})_{j \in \bb{\mathbf{s}}}$ in $\mathrm{Pth}_{\boldsymbol{\mathcal{A}}^{(2)}, \mathbf{s}}$ satisfying that, for evety $j \in \bb{\mathbf{s}}$,
$$
\mathrm{sc}^{(0,2)}_{\boldsymbol{\mathcal{A}}^{(2)}, s_{j}}\left(
\mathfrak{Q}^{(2)}_{j}
\right)
=
\mathrm{tg}^{(0,2)}_{\boldsymbol{\mathcal{A}}^{(2)}, s_{j}}\left(
\mathfrak{P}^{(2)}_{j}
\right);
$$
then
\begin{multline*}
\left(
\sigma^{\mathbf{Pth}_{\boldsymbol{\mathcal{A}}^{(2)}}}\left(
\left(
\mathfrak{Q}^{(2)}_{j} 
\circ_{s_{j}}^{0\mathbf{Pth}_{\boldsymbol{\mathcal{A}}^{(2)}}}
\mathfrak{P}^{(2)}_{j}
\right)_{j \in \bb{\mathbf{s}}}
\right)
,
\right.
\\
\left.
\sigma^{\mathbf{Pth}_{\boldsymbol{\mathcal{A}}^{(2)}}}\left(
\left(
\mathfrak{Q}^{(2)}_{j}
\right)_{j \in \bb{\mathbf{s}}}
\right)
\circ_{s}^{0\mathbf{Pth}_{\boldsymbol{\mathcal{A}}^{(2)}}}
\sigma^{\mathbf{Pth}_{\boldsymbol{\mathcal{A}}^{(2)}}}\left(
\left(
\mathfrak{P}^{(2)}_{j}
\right)_{j \in \bb{\mathbf{s}}}
\right)
\right)
\in
\Upsilon^{(1)}_{s}.
\end{multline*}

We want to check that
\begin{multline*}
\left\llbracket
f^{(2)\flat}_{s} \left(
\sigma^{\mathbf{Pth}_{\boldsymbol{\mathcal{A}}^{(2)}}}\left(
\left(
\mathfrak{Q}^{(2)}_{j} 
\circ_{s_{j}}^{0\mathbf{Pth}_{\boldsymbol{\mathcal{A}}^{(2)}}}
\mathfrak{P}^{(2)}_{j}
\right)_{j \in \bb{\mathbf{s}}}
\right)
\right)
\right\rrbracket_{\varphi(s)}
\\
=
\left\llbracket
f^{(2)\flat}_{s} \left(
\sigma^{\mathbf{Pth}_{\boldsymbol{\mathcal{A}}^{(2)}}}\left(
\left(
\mathfrak{Q}^{(2)}_{j}
\right)_{j \in \bb{\mathbf{s}}}
\right)
\circ_{s}^{0\mathbf{Pth}_{\boldsymbol{\mathcal{A}}^{(2)}}}
\sigma^{\mathbf{Pth}_{\boldsymbol{\mathcal{A}}^{(2)}}}\left(
\left(
\mathfrak{P}^{(2)}_{j}
\right)_{j \in \bb{\mathbf{s}}}
\right)
\right)
\right\rrbracket_{\varphi(s)}.
\end{multline*}

According to the fact that $\mathrm{pr}_{\boldsymbol{\mathcal{B}}^{(2)}, \varphi}^{\llbracket\cdot\rrbracket} \circ f^{(2)\flat}$ is a $\Sigma^{\boldsymbol{\mathcal{A}}^{(2)}}$-homomorphism, the desired equality is equivalent to check that
\begin{multline*}
\sigma^{\llbracket\mathbf{Pth}_{\boldsymbol{\mathcal{B}}^{(2)}}^{\mathbf{f}^{(2)}}\rrbracket}\left(
\left(
\left\llbracket
f^{(2)\flat}_{s_{j}} \left(
\mathfrak{Q}^{(2)}_{j}
\right)
\right\rrbracket_{\varphi(s_{j})}
\circ_{s_{j}}^{0\llbracket\mathbf{Pth}_{\boldsymbol{\mathcal{B}}^{(2)}}^{\mathbf{f}^{(2)}}\rrbracket}
\left\llbracket
f^{(2)\flat}_{s_{j}} \left(
\mathfrak{P}^{(2)}_{j}
\right)
\right\rrbracket_{\varphi(s_{j})}
\right)_{j \in \bb{\mathbf{s}}}
\right)
\\
=
\sigma^{\llbracket\mathbf{Pth}_{\boldsymbol{\mathcal{B}}^{(2)}}^{\mathbf{f}^{(2)}}\rrbracket}\left(
\left(
\left\llbracket
f^{(2)\flat}_{s_{j}} \left(
\mathfrak{Q}^{(2)}_{j}
\right)
\right\rrbracket_{\varphi(s_{j})}
\right)_{j \in \bb{\mathbf{s}}}
\right)
\circ_{s}^{0\llbracket\mathbf{Pth}_{\boldsymbol{\mathcal{B}}^{(2)}}^{\mathbf{f}^{(2)}}\rrbracket}
\\
\sigma^{\llbracket\mathbf{Pth}_{\boldsymbol{\mathcal{B}}^{(2)}}^{\mathbf{f}^{(2)}}\rrbracket}\left(
\left(
\left\llbracket
f^{(2)\flat}_{s_{j}} \left(
\mathfrak{P}^{(2)}_{j}
\right)
\right\rrbracket_{\varphi(s_{j})}
\right)_{j \in \bb{\mathbf{s}}}
\right)
.
\end{multline*}

But this equality follows from Proposition~\ref{PDQPthBVarA8}.

This proves Case (4).

This completes the proof of Proposition~\ref{PDHomPthExtKer}.
\end{proof}

Thus, according to the universal property of the quotient, we introduce the following definition.

\begin{definition}
\label{DDQPthExt}
Following Proposition~\ref{PDHomPthExtKer} and taking into account the Universal Property of the Quotient, there exists a unique $\Sigma^{\boldsymbol{\mathcal{A}}^{(2)}}$-homomorphism, that we will denote by $f^{\llbracket2\rrbracket @}$, i.e., 
$$
f^{\llbracket2\rrbracket @}
\colon
\llbracket \mathbf{Pth}_{\boldsymbol{\mathcal{A}^{(2)}}}\rrbracket
\mor
\llbracket \mathbf{Pth}_{\boldsymbol{\mathcal{B}^{(2)}}}^{\mathbf{f}^{(2)}}\rrbracket
$$
satisfying that $f^{\llbracket2\rrbracket @} \circ \mathrm{pr}^{\llbracket\cdot\rrbracket}_{\boldsymbol{\mathcal{A}}^{(2)}} = \mathrm{pr}_{\boldsymbol{\mathcal{B}}^{(2)}, \varphi}^{\llbracket\cdot\rrbracket} \circ f^{(2)\flat}$, namely $f^{\llbracket2\rrbracket @} = (\mathrm{pr}_{\boldsymbol{\mathcal{B}}^{(2)}, \varphi}^{\llbracket\cdot\rrbracket} \circ f^{(2)\flat})^{\natural}$. We will call this mapping the \emph{second-order quotient path extension mapping} of $\mathbf{f}^{(2)}$. Formally, for every sort $s$ in $S$ and every path class $\llbracket \mathfrak{P}^{(2)} \rrbracket_{s}$ in $\llbracket \mathrm{Pth}_{\boldsymbol{\mathcal{A}}^{(2)}} \rrbracket_{s}$, 
$$
f^{\llbracket2\rrbracket @}_{s}\left(
\left\llbracket
\mathfrak{P}^{(2)}
\right\rrbracket_{s}
\right)
=
\left\llbracket
f^{(2)\flat}_{s}\left(
\mathfrak{P}^{(2)}
\right)
\right\rrbracket_{\varphi(s)}
$$

Recall that, taking into account Theorems~\ref{TDIso} and \ref{TDIsoB},  we obtain a unique $\Sigma^{\boldsymbol{\mathcal{A}}^{(2)}}$-homomorphism from $\llbracket\mathbf{PT}_{\boldsymbol{\mathcal{A}}^{(2)}}\rrbracket$ to $\llbracket\mathbf{PT}_{\boldsymbol{\mathcal{B}}^{(2)}}^{\mathbf{f}^{(2)}}\rrbracket$. We agreed to also use $f^{\llbracket2\rrbracket @}$ to denote this mapping
$$
f^{\llbracket2\rrbracket @}
\colon
\llbracket \mathbf{PT}_{\boldsymbol{\mathcal{A}^{(2)}}}\rrbracket
\mor
\llbracket \mathbf{PT}_{\boldsymbol{\mathcal{B}^{(2)}}}^{\mathbf{f}^{(2)}}\rrbracket.
$$
Formally, for every sort $s$ in $S$ and every path term class $\llbracket P \rrbracket_{s}$ in $\llbracket \mathrm{PT}_{\boldsymbol{\mathcal{A}}^{(2)}}\rrbracket_{s}$,
$$
f^{\llbracket2\rrbracket @}_{s}\left(
\left\llbracket
P
\right\rrbracket_{s}
\right)
=
\left\llbracket
\mathrm{CH}^{(2)}_{\boldsymbol{\mathcal{B}}^{(2)}, \varphi(s)}\left(
f^{(2)\flat}_{s}\left(
\mathrm{ip}^{(2,X)@}_{\boldsymbol{\mathcal{A}}^{(2)}, s}\left(
P
\right)
\right)
\right)
\right\rrbracket_{\varphi(s)}.
$$

\end{definition}

\section{The behaviour of the second-order quotient path extension mapping}

To end this chapter, we study the relations between the second-order quotient path-extension mapping and the mappings $\mathrm{CH}^{\llbracket 2 \rrbracket}$ and $\mathrm{ip}^{(\llbracket 2 \rrbracket, X)@}$ and the source, target and identity path mappings.

\begin{proposition}
\label{PDQPthExtQCH}
Let $\mathbf{f}^{(2)}=(\varphi, c, (f^{(i)})_{i\in 3})$ be a morphism of first-order many-sorted rewriting systems from $\boldsymbol{\mathcal{A}}^{(2)}$ to $\boldsymbol{\mathcal{B}}^{(2)}$. Then the following equalities holds
$$
f^{\llbracket2\rrbracket @} \circ \mathrm{CH}_{\boldsymbol{\mathcal{A}}^{(2)}}^{\llbracket2\rrbracket}
=
\mathrm{CH}_{\boldsymbol{\mathcal{B}}^{(2)}, \varphi}^{\llbracket2\rrbracket} \circ f^{\llbracket2\rrbracket @}.
$$
\end{proposition}

\begin{proof}
For every sort $s$ in $S$ and every path class $\llbracket \mathfrak{P}^{(2)} \rrbracket_{s}$ in $\llbracket \mathrm{Pth}_{\boldsymbol{\mathcal{A}}^{(2)}} \rrbracket_{s}$, the following chain of equalities holds
\begin{flushleft}
$
f^{\llbracket 2 \rrbracket @}_{s} \left(
\mathrm{CH}^{\llbracket 2 \rrbracket}_{\boldsymbol{\mathcal{A}}^{(2)}, s} \left(
\left\llbracket
\mathfrak{P}^{(2)}
\right\rrbracket_{s}
\right)
\right)
$
\allowdisplaybreaks
\begin{align*}
&=
f^{\llbracket 2 \rrbracket @}_{s} \left(
\left\llbracket
\mathrm{CH}^{(2)}_{\boldsymbol{\mathcal{A}}^{(2)}, s} \left(
\mathfrak{P}^{(2)}
\right)
\right\rrbracket_{s}
\right)
\tag{1}
\\
&=
\left\llbracket
\mathrm{CH}^{(2)}_{\boldsymbol{\mathcal{B}}^{(2)}, \varphi(s)} \left(
f^{(2)\flat}_{s} \left(
\mathrm{ip}^{(2,X)@}_{\boldsymbol{\mathcal{A}}^{(2)}, s}\left(
\mathrm{CH}^{(2)}_{\boldsymbol{\mathcal{A}}^{(2)}, s} \left(
\mathfrak{P}^{(2)}
\right)
\right)
\right)
\right)
\right\rrbracket_{\varphi(s)}
\tag{2}
\\
&=
\mathrm{CH}^{\llbracket 2 \rrbracket}_{\boldsymbol{\mathcal{B}}^{(2)}, \varphi(s)} \left(
\left\llbracket
f^{\llbracket2\rrbracket @}_{s} \left(
\mathrm{ip}^{(2,X)@}_{\boldsymbol{\mathcal{A}}^{(2)}, s}\left(
\mathrm{CH}^{(2)}_{\boldsymbol{\mathcal{A}}^{(2)}, s} \left(
\mathfrak{P}^{(2)}
\right)
\right)
\right)
\right\rrbracket_{\varphi(s)}
\right)
\tag{3}
\\
&=
\mathrm{CH}^{\llbracket 2 \rrbracket}_{\boldsymbol{\mathcal{B}}^{(2)}, \varphi(s)} \left(
f^{\llbracket2\rrbracket @}_{s} \left(
\left\llbracket
\mathrm{ip}^{(2,X)@}_{\boldsymbol{\mathcal{A}}^{(2)}, s}\left(
\mathrm{CH}^{(2)}_{\boldsymbol{\mathcal{A}}^{(2)}, s} \left(
\mathfrak{P}^{(2)}
\right)
\right)
\right\rrbracket_{s}
\right)
\right)
\tag{4}
\\
&=
\mathrm{CH}^{\llbracket 2 \rrbracket}_{\boldsymbol{\mathcal{B}}^{(2)}, \varphi(s)} \left(
f^{\llbracket2\rrbracket @}_{s} \left(
\left\llbracket
\mathfrak{P}^{(2)}
\right\rrbracket_{s}
\right)
\right)
\tag{5}
\end{align*}
\end{flushleft}

The first equality unravels the definition of $s$-th component of the $S$-sorted mapping $\mathrm{CH}^{\llbracket 2 \rrbracket}_{\boldsymbol{\mathcal{A}}^{(2)}}$, introduced in Definition~\ref{DDPTQDCH};
the second equality unravels the definition of the $s$-th component of the $S$-sorted mapping $f^{\llbracket2\rrbracket @}$ at a path term class, introduced in Definition~\ref{DDQPthExt};
the third equality recovers the definition of the $\varphi(s)$-th component of the $S$-sorted mapping $\mathrm{CH}^{\llbracket 2 \rrbracket}_{\boldsymbol{\mathcal{B}}^{(2)}}$, introduced in Definition~\ref{DDPTQDCH};
the fourth equality recovers the definition of the $s$-th component of the $S$-sorted mapping $f^{\llbracket2\rrbracket @}$ at a path term class, introduced in Definition~\ref{DDQPthExt};
finally, the last equality follows from Proposition~\ref{PDIpDCH}.
\end{proof}

\begin{proposition}
\label{PDQPthExtQIp}
Let $\mathbf{f}^{(2)}=(\varphi, c, (f^{(i)})_{i\in 3})$ be a morphism of first-order many-sorted rewriting systems from $\boldsymbol{\mathcal{A}}^{(2)}$ to $\boldsymbol{\mathcal{B}}^{(2)}$. Then the following equalities holds
$$
\mathrm{ip}_{\boldsymbol{\mathcal{B}}^{(2)}, \varphi}^{(\llbracket2\rrbracket, Y)@}
\circ
f^{\llbracket2\rrbracket @}
=
f^{\llbracket2\rrbracket @}
\circ
\mathrm{ip}_{\boldsymbol{\mathcal{A}}^{(2)}}^{(\llbracket2\rrbracket, X)@}.
$$
\end{proposition}

\begin{proof}
For every sort $s$ in $S$ and every path term class $\llbracket P \rrbracket_{s}$ in $\llbracket \mathrm{PT}_{\boldsymbol{\mathcal{A}}^{(2)}} \rrbracket_{s}$, the following chain of equalities holds
\begin{flushleft}
$
\mathrm{ip}_{\boldsymbol{\mathcal{B}}^{(2)}, \varphi(s)}^{(\llbracket2\rrbracket, Y)@} \left(
f_{s}^{\llbracket2\rrbracket @} \left(
\left\llbracket
P
\right\rrbracket_{s}
\right)
\right)
$
\allowdisplaybreaks
\begin{align*}
&=
\mathrm{ip}_{\boldsymbol{\mathcal{B}}^{(2)}, \varphi(s)}^{(\llbracket2\rrbracket, Y)@} \left(
\left\llbracket
\mathrm{CH}^{(2)}_{\boldsymbol{\mathcal{B}}^{(2)}, \varphi(s)}\left(
f_{s}^{(2) \flat} \left(
\mathrm{ip}^{(2,X)@}_{\boldsymbol{\mathcal{A}}^{(2)}, s}\left(
P
\right)
\right)
\right)
\right\rrbracket_{\varphi(s)}
\right)
\tag{1}
\\
&=
\left\llbracket
\mathrm{ip}_{\boldsymbol{\mathcal{B}}^{(2)}, \varphi(s)}^{(2, Y)@} \left(
\mathrm{CH}^{(2)}_{\boldsymbol{\mathcal{B}}^{(2)}, \varphi(s)}\left(
f_{s}^{(2) \flat} \left(
\mathrm{ip}^{(2,X)@}_{\boldsymbol{\mathcal{A}}^{(2)}, s}\left(
P
\right)
\right)
\right)
\right)
\right\rrbracket_{\varphi(s)}
\tag{2}
\\
&=
\left\llbracket
f_{s}^{(2) \flat} \left(
\mathrm{ip}^{(2,X)@}_{\boldsymbol{\mathcal{A}}^{(2)}, s}\left(
P
\right)
\right)
\right\rrbracket_{\varphi(s)}
\tag{3}
\\
&=
f_{s}^{\llbracket 2 \rrbracket @} \left(
\left\llbracket
\mathrm{ip}^{(2,X)@}_{\boldsymbol{\mathcal{A}}^{(2)}, s}\left(
P
\right)
\right\rrbracket_{s}
\right)
\tag{4}
\\
&=
f_{s}^{\llbracket 2 \rrbracket @} \left(
\mathrm{ip}^{(2,X)@}_{\boldsymbol{\mathcal{A}}^{(2)}, s}\left(
\left\llbracket
P
\right\rrbracket_{s}
\right)
\right)
\tag{5}
\end{align*}
\end{flushleft}

The first equality unravels the definition of the $s$-th component of the $S$-sorted mapping $f^{\llbracket2\rrbracket @}$ at a path term class, introduced in Definition~\ref{DDQPthExt};
the second equality unravels the definition of the $\varphi(s)$-th component of the $S$-sorted mapping $\mathrm{ip}^{(\llbracket 2 \rrbracket,Y)@}_{\boldsymbol{\mathcal{B}}^{(2)}}$, introduced in Definition~\ref{DDPTQIp};
the third equality follows from Proposition~\ref{PDIpDCH};
the fourth equality recovers the definition of the $s$-th component of the $S$-sorted mapping $f^{\llbracket2\rrbracket @}$ at a path term class, introduced in Definition~\ref{DDQPthExt};
finally, the last equality recovers the definition of the $s$-th component of the $S$-sorted mapping $\mathrm{ip}^{(\llbracket 2 \rrbracket,X)@}_{\boldsymbol{\mathcal{A}}^{(2)}}$, introduced in Definition~\ref{DDPTQIp}.
\end{proof}

\begin{proposition}
\label{PDQPthExtDUScTg}
Let $\mathbf{f}^{(2)}=(\varphi, c, (f^{(i)})_{i\in 3})$ be a morphism of first-order many-sorted rewriting systems from $\boldsymbol{\mathcal{A}}^{(2)}$ to $\boldsymbol{\mathcal{B}}^{(2)}$. Then the following equalities holds
\begin{align*}
\mathrm{sc}_{\boldsymbol{\mathcal{B}}^{(2)},\varphi}^{([1],\llbracket2\rrbracket)} \circ f^{\llbracket2\rrbracket @}
&=
f^{[1] @}\circ\mathrm{sc}_{\boldsymbol{\mathcal{A}}^{(2)}}^{([1],\llbracket2\rrbracket)}
&&\mbox{and}&
\mathrm{tg}_{\boldsymbol{\mathcal{B}}^{(2)},\varphi}^{([1],\llbracket2\rrbracket)} \circ f^{\llbracket2\rrbracket @}
&=
f^{[1]@}\circ\mathrm{tg}_{\boldsymbol{\mathcal{A}}^{(2)}}^{([1],\llbracket2\rrbracket)}.
\end{align*}
\end{proposition}

\begin{proof}
For every sort $s$ in $S$ and every path class $\llbracket \mathfrak{P}^{(2)} \rrbracket_{s}$ in $\llbracket \mathrm{Pth}_{\boldsymbol{\mathcal{A}}^{(2)}} \rrbracket_{s}$, the following chain of equalities holds
\allowdisplaybreaks
\begin{align*}
\mathrm{sc}^{([1],\llbracket2\rrbracket)}_{\boldsymbol{\mathcal{B}}^{(2)}, \varphi(s)} \left(
f^{\llbracket 2 \rrbracket @}_{s} \left(
\left\llbracket \mathfrak{P}^{(2)} \right\rrbracket_{s}
\right)
\right)
&=
\mathrm{sc}^{([1],\llbracket2\rrbracket)}_{\boldsymbol{\mathcal{B}}^{(2)}, \varphi(s)} \left(
\left\llbracket
f^{(2)\flat}_{s} \left(
\mathfrak{P}^{(2)}
\right)
\right\rrbracket_{\varphi(s)}
\right)
\tag{1}
\\
&=
\mathrm{sc}^{([1],2)}_{\boldsymbol{\mathcal{B}}^{(2)}, \varphi(s)} \left(
f^{(2)\flat}_{s} \left(
\mathfrak{P}^{(2)}
\right)
\right)
\tag{2}
\\
&=
f^{[1]@}_{s} \left(
\mathrm{sc}^{([1],2)}_{\boldsymbol{\mathcal{A}}^{(2)}, s} \left(
\mathfrak{P}^{(2)}
\right)
\right)
\tag{3}
\\
&=
f^{[1]@}_{s} \left(
\mathrm{sc}^{([1],\llbracket2\rrbracket)}_{\boldsymbol{\mathcal{A}}^{(2)}, s} \left(
\left\llbracket
\mathfrak{P}^{(2)}
\right\rrbracket_{s}
\right)
\right)
\tag{4}
\end{align*}

The first equality unravels the definition of the second-order quotient path extension mapping $f^{\llbracket 2 \rrbracket @}$ at a path class, introduced in Definition~\ref{DDQPthExt};
the second equality unravels the definition of the mapping $\mathrm{sc}^{([1],\llbracket2\rrbracket)}_{\boldsymbol{\mathcal{B}}^{(2)}, \varphi}$, introduced in Definition~\ref{DDVDU};
the third equality follows from Proposition~\ref{PDPthExt};
finally, the last equality recovers the definition of the mapping $\mathrm{sc}^{([1],\llbracket2\rrbracket)}_{\boldsymbol{\mathcal{A}}^{(2)}}$, introduced in Definition~\ref{DDVDU}.

A similar argument applies to the $([1],\llbracket 2 \rrbracket)$-target mapping.
\end{proof}

\begin{proposition}
\label{PDQPthExtDUIp}
Let $\mathbf{f}^{(2)}=(\varphi, c, (f^{(i)})_{i\in 3})$ be a morphism of first-order many-sorted rewriting systems from $\boldsymbol{\mathcal{A}}^{(2)}$ to $\boldsymbol{\mathcal{B}}^{(2)}$. Then the following equality holds
$$
f^{\llbracket2\rrbracket @} \circ \mathrm{ip}_{\boldsymbol{\mathcal{A}}^{(2)}}^{(\llbracket2\rrbracket, [1])\sharp}
=
\mathrm{ip}_{\boldsymbol{\mathcal{B}}^{(2)}, \varphi}^{(\llbracket2\rrbracket,[1])\sharp} \circ f^{[1]@}.
$$
\end{proposition}

\begin{proof}
For every sort $s$ in $S$ and every path term class $[P]_{s}$ in $[\mathrm{PT}_{\boldsymbol{\mathcal{A}}^{(1)}}]_{s}$, the following chain of equalities holds
\allowdisplaybreaks
\begin{align*}
f^{\llbracket 2 \rrbracket @}_{s} \left(
\mathrm{ip}^{(\llbracket 2 \rrbracket, [1])\sharp}_{\boldsymbol{\mathcal{A}}^{(2)}, s}\left(
[P]_{s}
\right)
\right)
&=
f^{\llbracket 2 \rrbracket @}_{s} \left(
\left\llbracket
\mathrm{ip}^{(2, [1])\sharp}_{\boldsymbol{\mathcal{A}}^{(2)}, s}\left(
[P]_{s}
\right)
\right\rrbracket_{s}
\right)
\tag{1}
\\
&=
\left\llbracket
f^{(2)\flat}_{s} \left(
\mathrm{ip}^{(2, [1])\sharp}_{\boldsymbol{\mathcal{A}}^{(2)}, s}\left(
[P]_{s}
\right)
\right)
\right\rrbracket_{\varphi(s)}
\tag{2}
\\
&=
\left\llbracket
\mathrm{ip}^{(2, [1])\sharp}_{\boldsymbol{\mathcal{B}}^{(2)}, \varphi(s)}\left(
f^{[1]@}_{s} \left(
[P]_{s}
\right)
\right)
\right\rrbracket_{\varphi(s)}
\tag{3}
\\
&=
\mathrm{ip}^{(\llbracket 2 \rrbracket, [1])\sharp}_{\boldsymbol{\mathcal{B}}^{(2)}, \varphi(s)}\left(
f^{[1]@}_{s} \left(
[P]_{s}
\right)
\right)
\tag{4}
\end{align*}

The first equality unravels the definition of the mapping $\mathrm{ip}^{(\llbracket2\rrbracket, [1])}_{\boldsymbol{\mathcal{A}}^{(2)}}$, introduced in Definition~\ref{DDVDU};
the second equality unravels the definition of the second-order quotient path extension mapping $f^{\llbracket 2 \rrbracket @}$ at a path class, introduced in Definition~\ref{DDQPthExt};
the third equality follows from Proposition~\ref{PDPthExt};
finally, the last equality recovers the definition of the mapping $\mathrm{ip}^{(\llbracket2\rrbracket, [1])}_{\boldsymbol{\mathcal{B}}^{(2)}, \varphi}$, introduced in Definition~\ref{DDVDU}.
\end{proof}

\begin{proposition}
\label{PDQPthExtDZScTg}
Let $\mathbf{f}^{(2)}=(\varphi, c, (f^{(i)})_{i\in 3})$ be a morphism of first-order many-sorted rewriting systems from $\boldsymbol{\mathcal{A}}^{(2)}$ to $\boldsymbol{\mathcal{B}}^{(2)}$. Then the following equalities holds
\begin{align*}
\mathrm{sc}_{\boldsymbol{\mathcal{B}}^{(2)},\varphi}^{(0,\llbracket2\rrbracket)} \circ f^{\llbracket2\rrbracket @}
&=
f^{(0)\sharp}\circ\mathrm{sc}_{\boldsymbol{\mathcal{A}}^{(2)}}^{(0,\llbracket2\rrbracket)}
&&\mbox{and}&
\mathrm{tg}_{\boldsymbol{\mathcal{B}}^{(2)},\varphi}^{(0,\llbracket2\rrbracket)} \circ f^{\llbracket2\rrbracket @}
&=
f^{(0)\sharp}\circ\mathrm{tg}_{\boldsymbol{\mathcal{A}}^{(2)}}^{(0,\llbracket2\rrbracket)}.
\end{align*}
\end{proposition}

\begin{proof}
For every sort $s$ in $S$ and every path class $\llbracket \mathfrak{P}^{(2)} \rrbracket_{s}$ in $\llbracket \mathrm{Pth}_{\boldsymbol{\mathcal{A}}^{(2)}} \rrbracket_{s}$, the following chain of equalities holds
\allowdisplaybreaks
\begin{align*}
\mathrm{sc}^{(0,\llbracket2\rrbracket)}_{\boldsymbol{\mathcal{B}}^{(2)}, \varphi(s)} \left(
f^{\llbracket 2 \rrbracket @}_{s} \left(
\left\llbracket \mathfrak{P}^{(2)} \right\rrbracket_{s}
\right)
\right)
&=
\mathrm{sc}^{(0,\llbracket2\rrbracket)}_{\boldsymbol{\mathcal{B}}^{(2)}, \varphi(s)} \left(
\left\llbracket
f^{(2)\flat}_{s} \left(
\mathfrak{P}^{(2)}
\right)
\right\rrbracket_{\varphi(s)}
\right)
\tag{1}
\\
&=
\mathrm{sc}^{(0,2)}_{\boldsymbol{\mathcal{B}}^{(2)}, \varphi(s)} \left(
f^{(2)\flat}_{s} \left(
\mathfrak{P}^{(2)}
\right)
\right)
\tag{2}
\\
&=
f^{(0)\sharp}_{s} \left(
\mathrm{sc}^{(0,2)}_{\boldsymbol{\mathcal{A}}^{(2)}, s} \left(
\mathfrak{P}^{(2)}
\right)
\right)
\tag{3}
\\
&=
f^{(0)\sharp}_{s} \left(
\mathrm{sc}^{(0,\llbracket2\rrbracket)}_{\boldsymbol{\mathcal{A}}^{(2)}, s} \left(
\left\llbracket
\mathfrak{P}^{(2)}
\right\rrbracket_{s}
\right)
\right)
\tag{4}
\end{align*}

The first equality unravels the definition of the second-order quotient path extension mapping $f^{\llbracket 2 \rrbracket @}$ at a path class, introduced in Definition~\ref{DDQPthExt};
the second equality unravels the definition of the mapping $\mathrm{sc}^{(0,\llbracket2\rrbracket)}_{\boldsymbol{\mathcal{B}}^{(2)}, \varphi}$, introduced in Definition~\ref{DDVDZ};
the third equality follows from Proposition~\ref{PDPthExtDZScTg};
finally, the last equality recovers the definition of the mapping $\mathrm{sc}^{(0,\llbracket2\rrbracket)}_{\boldsymbol{\mathcal{A}}^{(2)}}$, introduced in Definition~\ref{DDVDZ}.

A similar argument applies to the $(0,\llbracket 2 \rrbracket)$-target mapping.
\end{proof}

\begin{proposition}
\label{PDQPthExtDZIp}
Let $\mathbf{f}^{(2)}=(\varphi, c, (f^{(i)})_{i\in 3})$ be a morphism of first-order many-sorted rewriting systems from $\boldsymbol{\mathcal{A}}^{(2)}$ to $\boldsymbol{\mathcal{B}}^{(2)}$. Then the following equality holds
$$
f^{\llbracket2\rrbracket @} \circ \mathrm{ip}_{\boldsymbol{\mathcal{A}}^{(2)}}^{(\llbracket2\rrbracket, 0)\sharp}
=
\mathrm{ip}_{\boldsymbol{\mathcal{B}}^{(2)}, \varphi}^{(\llbracket2\rrbracket, 0)\sharp} \circ f^{(0)\sharp}.
$$
\end{proposition}

\begin{proof}
For every sort $s$ in $S$ and every term $P$ in $\T_{\Sigma}(X)$, the following chain of equalities holds
\allowdisplaybreaks
\begin{align*}
f^{\llbracket 2 \rrbracket @}_{s} \left(
\mathrm{ip}^{(\llbracket 2 \rrbracket, 0)\sharp}_{\boldsymbol{\mathcal{A}}^{(2)}, s}\left(
P
\right)
\right)
&=
f^{\llbracket 2 \rrbracket @}_{s} \left(
\left\llbracket
\mathrm{ip}^{(2, 0)\sharp}_{\boldsymbol{\mathcal{A}}^{(2)}, s}\left(
P
\right)
\right\rrbracket_{s}
\right)
\tag{1}
\\
&=
\left\llbracket
f^{(2)\flat}_{s} \left(
\mathrm{ip}^{(2, 0)\sharp}_{\boldsymbol{\mathcal{A}}^{(2)}, s}\left(
P
\right)
\right)
\right\rrbracket_{\varphi(s)}
\tag{2}
\\
&=
\left\llbracket
\mathrm{ip}^{(2, 0)\sharp}_{\boldsymbol{\mathcal{B}}^{(2)}, \varphi(s)}\left(
f^{[1]@}_{s} \left(
P
\right)
\right)
\right\rrbracket_{\varphi(s)}
\tag{3}
\\
&=
\mathrm{ip}^{(\llbracket 2 \rrbracket, 0)\sharp}_{\boldsymbol{\mathcal{B}}^{(2)}, \varphi(s)}\left(
f^{[1]@}_{s} \left(
P
\right)
\right)
\tag{4}
\end{align*}

The first equality unravels the definition of the mapping $\mathrm{ip}^{(\llbracket2\rrbracket, 0)}_{\boldsymbol{\mathcal{A}}^{(2)}}$, introduced in Definition~\ref{DDVDZ};
the second equality unravels the definition of the second-order quotient path extension mapping $f^{\llbracket 2 \rrbracket @}$ at a path class, introduced in Definition~\ref{DDQPthExt};
the third equality follows from Proposition~\ref{PDPthExtDZIp};
finally, the last equality recovers the definition of the mapping $\mathrm{ip}^{(\llbracket2\rrbracket, 0)}_{\boldsymbol{\mathcal{B}}^{(2)}, \varphi}$, introduced in Definition~\ref{DDVDZ}.
\end{proof}			
\chapter{The category of second-order many-sorted rewriting systems}\label{S3H}

In this chapter we introduce the notions of second-order identity morphism at a second-order many-sorted rewriting system and that of second-order composition morphism of second-order morphisms. We then define the second-order equivalence relation between second-order morphisms to show that second-order many-sorted rewriting systems with equivalence classes of second-order morphisms form a category denoted $\mathsf{Rws}_{\mathfrak{d}}^{\llbracket 2 \rrbracket}$. Moreover, we introduce the notion of second-order tower associated to a second-order rewriting system and that of morphism of second-order towers. Then showing that the second-order towers and morphisms between them form a category denoted $\mathsf{Tw}_{\mathfrak{d}}^{\llbracket 2 \rrbracket}$. Indeed we show that the categories $\mathsf{Rws}_{\mathfrak{d}}^{\llbracket 2 \rrbracket}$ and $\mathsf{Tw}_{\mathfrak{d}}^{\llbracket 2 \rrbracket}$ are isomorphic. 

\section{The category $\mathsf{Rws}_{\mathfrak{d}}^{\llbracket 2 \rrbracket}$}

We begin by defining the identity morphism at a secod-order rewriting system and studying its properties.

\begin{definition}
\label{DIdRws2}
Let $\boldsymbol{\mathcal{A}}^{(2)} = (\boldsymbol{\mathcal{A}}^{(1)}, \mathcal{A}^{(2)})$ be a second-order many-sorted rewriting system. The \emph{second-order identity morphism} at $\boldsymbol{\mathcal{A}}^{(2)}$, denoted by $\mathrm{id}^{\boldsymbol{\mathcal{A}}^{(2)}}$, is given by 
$$
\mathrm{id}^{\boldsymbol{\mathcal{A}}^{(2)}}
=
\left(
\mathrm{id}^{\boldsymbol{\mathcal{A}}^{(1)}}, 
\mathrm{ech}^{(2,\mathcal{A}^{(2)})}_{\boldsymbol{\mathcal{A}}^{(2)}}
\right).
$$
where $\mathrm{id}^{\boldsymbol{\mathcal{A}}^{(1)}}$ is the identity first-order morphism at the first-order many-sorted rewriting system $\boldsymbol{\mathcal{A}}^{(1)}$, introduced in Definition~\ref{DIdRws1}, and $\mathrm{ech}^{(2,\mathcal{A}^{(2)})}_{\boldsymbol{\mathcal{A}}^{(2)}}$ is the $S$-sorted second-order echelon mapping associated with the many-sorted rewriting system $\boldsymbol{\mathcal{A}}^{(2)}$, introduced in Definition~\ref{DDEch}.
\end{definition}

\begin{proposition}
\label{PDPthExtEch}
Let $\boldsymbol{\mathcal{A}}^{(2)}$ be a second-order rewriting system and let $\mathrm{id}^{\boldsymbol{\mathcal{A}}^{(2)}}$ be its second-order identity morphism. Thus, $\mathrm{ech}^{(2,\mathcal{A}^{(2)})\flat}_{\boldsymbol{\mathcal{A}}^{(2)}} = \mathrm{ip}_{\boldsymbol{\mathcal{A}}^{(2)}}^{(2,X)@} \circ \mathrm{CH}^{(2)}_{\boldsymbol{\mathcal{A}}^{(2)}}$.
\end{proposition}

\begin{proof}
We want to show that, for every sort $s \in S$ and every second-order path $\mathfrak{P}^{(2)}$ in $\mathrm{Pth}_{\boldsymbol{\mathcal{A}}^{(2)},s}$, the path $\mathrm{ech}^{(2,\mathcal{A}^{(2)})\flat}_{\boldsymbol{\mathcal{A}}^{(2)}} (\mathfrak{P}^{(2)})$ is equal to $\mathrm{ip}_{\boldsymbol{\mathcal{A}}^{(2)}}^{(2,X)@} \circ \mathrm{CH}^{(2)}_{\boldsymbol{\mathcal{A}}^{(2)}} (\mathfrak{P}^{(2)})$.

We prove the statement by Artinian induction on $(\coprod\mathrm{Pth}_{\boldsymbol{\mathcal{A}}^{(2)}}, \leq_{\mathbf{Pth}_{\boldsymbol{\mathcal{A}}^{(2)}}})$.

\textsf{Base step of the Artinian induction.}

Let $(\mathfrak{P}^{(2)}, s)$ be a minimal element in $(\coprod\mathrm{Pth}_{\boldsymbol{\mathcal{A}}^{(2)}}, \leq_{\mathbf{Pth}_{\boldsymbol{\mathcal{A}}^{(2)}}})$. Then by Proposition~\ref{PDMinimal}, it follows that $\mathfrak{P}^{(2)}$ is either (1) an $(2,[1])$-identity second-order path or (2) a second-order echelon.

If~(1), i.e., if $\mathfrak{P}^{(2)}$ is a $(2,[1])$-identity path, then $\mathfrak{P}^{(2)} = \mathrm{ip}^{(2,[1])\sharp}_{\boldsymbol{\mathcal{A}}^{(2)}, s}([P]_{s})$ for some path term class $[P]_{s} \in [\mathrm{PT}_{\boldsymbol{\mathcal{A}}^{(2)}}]_{s}$. Thus, the following chain of equalities holds
\allowdisplaybreaks
\begin{align*}
\mathrm{ech}^{(2,\mathcal{A}^{(2)})\flat}_{\boldsymbol{\mathcal{A}}^{(2)}, s}(\mathfrak{P}^{(2)})
&=
\mathrm{ech}^{(2,\mathcal{A}^{(2)})\flat}_{\boldsymbol{\mathcal{A}}^{(2)}, s}\left(
\mathrm{ip}^{(2,[1])\sharp}_{\boldsymbol{\mathcal{A}}^{(2)}, s}\left(
[P]_{2}
\right)
\right)
\tag{1}
\\
&=
\mathrm{ip}^{(2,[1])\sharp}_{\boldsymbol{\mathcal{A}}^{(2)}, s}\left(
\mathrm{ech}^{(1, \mathcal{A}^{(1)})@}_{\boldsymbol{\mathcal{A}}^{(1)}}\left(
[P]_{s}
\right)
\right)
\tag{2}
\\
&=
\mathrm{ip}^{(2,[1])\sharp}_{\boldsymbol{\mathcal{A}}^{(2)}, s}\left(
[P]_{s}
\right)
\tag{3}
\\
&=
\eta^{(2,\mathbf{Pth}_{\boldsymbol{\mathcal{A}}^{(2)}})}\left(
\mathrm{ip}^{(2,[1])\sharp}_{\boldsymbol{\mathcal{A}}^{(2)}, s}\left(
[P]_{s}
\right)
\right)
\tag{4}
\\
&=
\mathrm{ip}^{(2,X)@}_{\boldsymbol{\mathcal{A}}^{(2)}, s}\left(
\mathrm{CH}^{(2)}_{\boldsymbol{\mathcal{A}}^{(2)}, s}\left(
\mathrm{ip}^{(2,[1])\sharp}_{\boldsymbol{\mathcal{A}}^{(2)}, s}\left(
[P]_{s}
\right)
\right)
\right)
\tag{5}
\\
&=
\mathrm{ip}^{(2,X)@}_{\boldsymbol{\mathcal{A}}^{(2)}, s}\left(
\mathrm{CH}^{(2)}_{\boldsymbol{\mathcal{A}}^{(2)}, s}\left(
\mathfrak{P}^{(2)}
\right)
\right).
\tag{6}
\end{align*}

In the just stated chain of equalities, the first equality unravels the definition of the path $\mathfrak{P}^{(2)}$;
the second equality unravels the definition of the second-order path extension mapping at a $(2,[1])$-identity path, introduced in Proposition~\ref{PDPthExt};
the third equality follows from Proposition~\ref{PQPthExtEch};
the fourth equality follows from the fact that we identify $\eta^{(2,\mathbf{Pth}_{\boldsymbol{\mathcal{A}}^{(2)}})}(\mathfrak{P}^{(2)})$ with $\mathfrak{P}^{(2)}$;
the fifth equality follows from Proposition~\ref{PDIpDU};
finally, the last equality recovers the definition of the path $\mathfrak{P}^{(2)}$.

If~(2), i.e., if $\mathfrak{P}^{(2)}$ is a second-order echelon associated to a second-order rewrite rule $\mathfrak{p}^{(2)}$ in $\mathcal{A}^{(2)}_{s}$, then the following chain of equalities holds
\allowdisplaybreaks
\begin{align*}
\mathrm{ech}^{(2,\mathcal{A}^{(2)})\flat}_{\boldsymbol{\mathcal{A}}^{(2)}, s}(\mathfrak{P}^{(2)})
&=
\mathrm{ech}^{(2,\mathcal{A}^{(2)})\flat}_{\boldsymbol{\mathcal{A}}^{(2)}, s}\left(
\mathrm{ech}^{(2,\mathcal{A}^{(2)})}_{\boldsymbol{\mathcal{A}}^{(2)}, s}\left(\mathfrak{p}^{(2)}\right)
\right)
\tag{1}
\\
&=
\mathrm{ech}^{(2,\mathcal{A}^{(2)})}_{\boldsymbol{\mathcal{A}}^{(2)}, s}\left(
\mathfrak{p}^{(2)}
\right)
\tag{2}
\\
&=
\mathfrak{p}^{(2)\mathbf{Pth}_{\boldsymbol{\mathcal{A}}^{(2)}}}
\tag{3}
\\
&=
\mathfrak{p}^{(2)\mathbf{F}_{\Sigma^{\boldsymbol{\mathcal{A}}^{(2)}}}(\mathbf{Pth}_{\boldsymbol{\mathcal{A}}^{(2)}})}
\tag{4}
\\
&=
\mathrm{ip}^{(2,X)@}_{\boldsymbol{\mathcal{A}}^{(2)}, s}\left(
\mathfrak{p}^{(2)\mathbf{T}_{\Sigma^{\boldsymbol{\mathcal{A}}^{(2)}}}(X)}
\right)
\tag{5}
\\
&=
\mathrm{ip}^{(2,X)@}_{\boldsymbol{\mathcal{A}}^{(2)}, s}\left(
\mathrm{CH}^{(2)}_{\boldsymbol{\mathcal{A}}^{(2)}, s} \left(
\mathrm{ech}^{(2,\mathcal{A}^{(2)})}_{\boldsymbol{\mathcal{A}}^{(2)}, s}\left(\mathfrak{p}^{(2)}\right)
\right)
\right)
\tag{6}
\\
&=
\mathrm{ip}^{(2,X)@}_{\boldsymbol{\mathcal{A}}^{(2)}, s}\left(
\mathrm{CH}^{(2)}_{\boldsymbol{\mathcal{A}}^{(2)}, s} \left(
\mathfrak{P}^{(2)}
\right)
\right).
\tag{7}
\end{align*}

In the just stated chain of equalities, the first equality unravels the definition of $\mathfrak{P}^{(2)}$;
the second equality unravels the definition of the second-order path extension mapping at a second-order echelon, introduced in Proposition~\ref{PDPthExt};
the third equality recovers the interpretation of the constant symbol $\mathfrak{p}^{(2)}$ in the many-sorted partial $\Sigma^{\boldsymbol{\mathcal{A}}^{(2)}}$-algebra introduced in Proposition~\ref{PDPthDCatAlg};
the fourth equality holds since, by Proposition~\ref{PDPthDCatAlg}, we have an interpretation of the constant symbol $\mathfrak{p}^{(2)}$;
the fifth equality holds because, according to Definition~\ref{DDIp}, $\mathrm{ip}^{(2,X)@}_{\boldsymbol{\mathcal{A}}^{(2)}}$ is a many-sorted $\Sigma^{\boldsymbol{\mathcal{A}}^{(2)}}$-homomorphism;
the sixth equality recovers the definition of the second-order Curry-Howard mapping at a second-order echelon, introduced in Definition~\ref{DDCH};
finally, the last equality recovers the definition of $\mathfrak{P}^{(2)}$.

This proves the base case of the induction.

\textsf{Inductive step of the Artinian induction.}

Let $(\mathfrak{P}^{(2)},s)$ be a non-minimal element of $(\coprod\mathrm{Pth}_{\boldsymbol{\mathcal{A}}^{(2)}}, \leq_{\mathbf{Pth}_{\boldsymbol{\mathcal{A}}^{(2)}}})$. We can assume that $\mathfrak{P}^{(2)}$ is not a $(2,[1])$-identity second-order path, since that case has already been considered. Let us suppose that, for every sort $t \in S$ and every second-order path $\mathfrak{Q}^{(2)} \in \mathrm{Pth}_{\boldsymbol{\mathcal{A}}^{(2)}, t}$, if $(\mathfrak{Q}^{(2)}, t) <_{\mathbf{Pth}_{\boldsymbol{\mathcal{A}}^{(2)}}} (\mathfrak{P}^{(2)}, s)$, then the value of $\mathrm{ech}^{(2,\mathcal{A}^{(2)})\flat}_{\boldsymbol{\mathcal{A}}^{(2)}, t}(\mathfrak{Q}^{(2)})$ is equal to $\mathrm{ip}^{(2,X)@}_{\boldsymbol{\mathcal{A}}^{(2)}, t}(
\mathrm{CH}^{(2)}_{\boldsymbol{\mathcal{A}}^{(2)}, t}(
\mathfrak{Q}^{(2)}
)
)$. 

By Lemma~\ref{LDOrdI}, we have that $\mathfrak{P}^{(2)}$ is either~(1) a second-order path of length strictly greater that one containing at least one second-order echelon or~(2) an echelonless second-order path.

If~(1), i.e., if $\mathfrak{P}^{(2)}$ is a second-order path of legth strictly greater than one containing at least one second-order echelon. Then, we let $i \in \bb{\mathfrak{P}^{(2)}}$ be the first index for which the one-step subpath $\mathfrak{P}^{(2),i,i}$ is a second-order echelon. We distinguish two cases accordingly, either~(1.1) $i=0$ or~(1.2) $i>0$. 

If~(1.1), i.e., if $\mathfrak{P}^{(2)}$ is a path of length strictly greater having its first second-order echelon on its first step, then the following chain of equalities holds
\begin{flushleft}
$
\mathrm{ech}^{(2,\mathcal{A}^{(2)})\flat}_{\boldsymbol{\mathcal{A}}^{(2)}, s}(\mathfrak{P}^{(2)})
$
\allowdisplaybreaks
\begin{align*}
&=
\mathrm{ech}^{(2,\mathcal{A}^{(2)})\flat}_{\boldsymbol{\mathcal{A}}^{(2)}, s}\left(
\mathfrak{P}^{(2), 1, \bb{\mathfrak{P}^{(2)}}-1}
\right)
\circ^{1\mathbf{Pth}_{\boldsymbol{\mathcal{A}}^{(2)}}}_{s}
\mathrm{ech}^{(2,\mathcal{A}^{(2)})\flat}_{\boldsymbol{\mathcal{A}}^{(2)}, s}\left(
\mathfrak{P}^{(2),0,0}
\right)
\tag{1}
\\
&=
\resizebox{.89\textwidth}{!}{$
\mathrm{ip}^{(2,X)@}_{\boldsymbol{\mathcal{A}}^{(2)}, s}\left(
\mathrm{CH}^{(2)}_{\boldsymbol{\mathcal{A}}^{(2)}, s} \left(
\mathfrak{P}^{(2), 1, \bb{\mathfrak{P}^{(2)}}-1}
\right)
\right)
\circ^{1\mathbf{Pth}_{\boldsymbol{\mathcal{A}}^{(2)}}}_{s}
\mathrm{ip}^{(2,X)@}_{\boldsymbol{\mathcal{A}}^{(2)}, s}\left(
\mathrm{CH}^{(2)}_{\boldsymbol{\mathcal{A}}^{(2)}, s} \left(
\mathfrak{P}^{(2),0,0}
\right)
\right)
$}
\tag{2}
\\
&=
\resizebox{.89\textwidth}{!}{$
\mathrm{ip}^{(2,X)@}_{\boldsymbol{\mathcal{A}}^{(2)}, s}\left(
\mathrm{CH}^{(2)}_{\boldsymbol{\mathcal{A}}^{(2)}, s} \left(
\mathfrak{P}^{(2), 1, \bb{\mathfrak{P}^{(2)}}-1}
\right)
\right)
\circ^{1\mathbf{F}_{\Sigma^{\boldsymbol{\mathcal{A}}^{(2)}}}(\mathbf{Pth}_{\boldsymbol{\mathcal{A}}^{(2)}})}_{s}
\mathrm{ip}^{(2,X)@}_{\boldsymbol{\mathcal{A}}^{(2)}, s}\left(
\mathrm{CH}^{(2)}_{\boldsymbol{\mathcal{A}}^{(2)}, s} \left(
\mathfrak{P}^{(2),0,0}
\right)
\right)
$}
\tag{3}
\\
&=
\mathrm{ip}^{(2,X)@}_{\boldsymbol{\mathcal{A}}^{(2)}, s}\left(
\mathrm{CH}^{(2)}_{\boldsymbol{\mathcal{A}}^{(2)}, s} \left(
\mathfrak{P}^{(2), 1, \bb{\mathfrak{P}^{(2)}}-1}
\right)
\circ^{1\mathbf{T}_{\Sigma^{\boldsymbol{\mathcal{A}}^{(2)}}}(X)}_{s}
\mathrm{CH}^{(2)}_{\boldsymbol{\mathcal{A}}^{(2)}, s} \left(
\mathfrak{P}^{(2),0,0}
\right)
\right)
\tag{4}
\\
&=
\mathrm{ip}^{(2,X)@}_{\boldsymbol{\mathcal{A}}^{(2)}, s}\left(
\mathrm{CH}^{(2)}_{\boldsymbol{\mathcal{A}}^{(2)}, s} \left(
\mathfrak{P}^{(2)}
\right)
\right).
\tag{5}
\end{align*}
\end{flushleft}

In the just stated chain of equalities, the first equality unravels the definition of the second-order path extension mapping at a path of length strictly greater than one with a second-order echelon on its first step, introduced in Proposition~\ref{PDPthExt};
the second equality follows taking into account that $(\mathfrak{P}^{((2),1,\bb{\mathfrak{P}^{(2)}}-1}, s)$ and $(\mathfrak{P}^{(2),0,0}, s)$ are strictly smaller that $(\mathfrak{P}^{(2)}, s)$ with respect to $\prec_{\mathbf{Pth}_{\boldsymbol{\mathcal{A}}^{(2)}}}$;
the third equality holds because, according to Proposition~\ref{PDIpDCH}, we have that 
$$
\mathrm{ip}^{(2,X)@}_{\boldsymbol{\mathcal{A}}^{(2)}, s}(
\mathrm{CH}^{(2)}_{\boldsymbol{\mathcal{A}}^{(2)}, s} (
\mathfrak{P}^{(2), 1, \bb{\mathfrak{P}^{(2)}}-1}
)
)
\mbox{ and }
\mathrm{ip}^{(2,X)@}_{\boldsymbol{\mathcal{A}}^{(2)}, s}(
\mathrm{CH}^{(2)}_{\boldsymbol{\mathcal{A}}^{(2)}, s}(
\mathfrak{P}^{(2),0,0}
)
)
$$
are paths in $\llbracket\mathfrak{P}^{(2),1, \bb{\mathfrak{P}^{(2)}}-1}\rrbracket_{s}$ and $\llbracket\mathfrak{P}^{(2),0,0}\rrbracket_{s}$, respectively. Hence the interpretation of the $1$-composition operation symbol $\circ^{1}_{s}$ in $\mathbf{F}_{\Sigma^{\boldsymbol{\mathcal{A}}^{(2)}}}(\mathbf{Pth}_{\boldsymbol{\mathcal{A}}^{(2)}})$ becomes that of $\mathbf{Pth}_{\boldsymbol{\mathcal{A}}^{(2)}}$.
The fourth equality holds since, according to Definition~\ref{DDIp}, $\mathrm{ip}^{(2,X)@}_{\boldsymbol{\mathcal{A}}^{(2)}}$ is a many-sorted $\Sigma^{\boldsymbol{\mathcal{A}}^{(2)}}$-homomorphism;
finally, the last equality recovers the definition of the second-order Curry-Howard mapping at a path of length strictly greater than one with a second-order echelon on its first step, introduced in Definition~\ref{DDCH}.

If~(1.2), i.e., if $\mathfrak{P}^{(2)}$ is a second-order path of length strictly greater having its first second-order echelon at position $i \in \bb{\mathfrak{P}^{(2)}}$, with $i > 0$, then the following chain of equalities holds
\begin{flushleft}
$
\mathrm{ech}^{(2,\mathcal{A}^{(2)})\flat}_{\boldsymbol{\mathcal{A}}^{(2)}, s}(\mathfrak{P}^{(2)})
$
\allowdisplaybreaks
\begin{align*}
&=
\mathrm{ech}^{(2,\mathcal{A}^{(2)})\flat}_{\boldsymbol{\mathcal{A}}^{(2)}, s}\left(
\mathfrak{P}^{(2),i, \bb{\mathfrak{P}^{(2)}}-1}
\right)
\circ^{1\mathbf{Pth}_{\boldsymbol{\mathcal{A}}^{(1)}}}_{s}
\mathrm{ech}^{(2,\mathcal{A}^{(2)})\flat}_{\boldsymbol{\mathcal{A}}^{(2)}, s}\left(
\mathfrak{P}^{(2),0,i-1}
\right)
\tag{1}
\\
&=
\resizebox{.89\textwidth}{!}{$
\mathrm{ip}^{(2,X)@}_{\boldsymbol{\mathcal{A}}^{(2)}, s}\left(
\mathrm{CH}^{(2)}_{\boldsymbol{\mathcal{A}}^{(2)}, s} \left(
\mathfrak{P}^{(2), i, \bb{\mathfrak{P}^{(2)}}-1}
\right)
\right)
\circ^{1\mathbf{Pth}_{\boldsymbol{\mathcal{A}}^{(1)}}}_{s}
\mathrm{ip}^{(2,X)@}_{\boldsymbol{\mathcal{A}}^{(2)}, s}\left(
\mathrm{CH}^{(2)}_{\boldsymbol{\mathcal{A}}^{(2)}, s} \left(
\mathfrak{P}^{(2),0,i-1}
\right)
\right)
$}
\tag{2}
\\
&=
\resizebox{.89\textwidth}{!}{$
\mathrm{ip}^{(2,X)@}_{\boldsymbol{\mathcal{A}}^{(2)}, s}\left(
\mathrm{CH}^{(2)}_{\boldsymbol{\mathcal{A}}^{(2)}, s} \left(
\mathfrak{P}^{(2), i, \bb{\mathfrak{P}^{(2)}}-1}
\right)
\right)
\circ^{1\mathbf{F}_{\Sigma^{\boldsymbol{\mathcal{A}}^{(2)}}}(\mathbf{Pth}_{\boldsymbol{\mathcal{A}}^{(2)}})}_{s}
\mathrm{ip}^{(2,X)@}_{\boldsymbol{\mathcal{A}}^{(2)}, s}\left(
\mathrm{CH}^{(2)}_{\boldsymbol{\mathcal{A}}^{(2)}, s} \left(
\mathfrak{P}^{(2),0,i-1}
\right)
\right)
$}
\tag{3}
\\
&=
\mathrm{ip}^{(2,X)@}_{\boldsymbol{\mathcal{A}}^{(2)}, s}\left(
\mathrm{CH}^{(2)}_{\boldsymbol{\mathcal{A}}^{(2)}, s} \left(
\mathfrak{P}^{(2),i, \bb{\mathfrak{P}^{(2)}}-1}
\right)
\circ^{1\mathbf{T}_{\Sigma^{\boldsymbol{\mathcal{A}}^{(2)}}}(X)}_{s}
\mathrm{CH}^{(2)}_{\boldsymbol{\mathcal{A}}^{(2)}, s} \left(
\mathfrak{P}^{(2),0,i-1}
\right)
\right)
\tag{4}
\\
&=
\mathrm{ip}^{(2,X)@}_{\boldsymbol{\mathcal{A}}^{(2)}, s}\left(
\mathrm{CH}^{(2)}_{\boldsymbol{\mathcal{A}}^{(2)}, s} \left(
\mathfrak{P}^{(2)}
\right)
\right).
\tag{5}
\end{align*}
\end{flushleft}

In the just stated chain of equalities, the first equality unravels the definition of the second-order path extension mapping at a second-order path of length strictly greater than one with a second-order echelon at position $i$, introduced in Proposition~\ref{PDPthExt};
the second equality follows taking into account that $(\mathfrak{P}^{((2),i,\bb{\mathfrak{P}^{(2)}}-1}, s)$ and $(\mathfrak{P}^{(2),0,i-1}, s)$ are strictly smaller that $(\mathfrak{P}^{(2)}, s)$ with respect to $\prec_{\mathbf{Pth}_{\boldsymbol{\mathcal{A}}^{(2)}}}$;
the third equality holds because, according to Proposition~\ref{PDIpDCH}, we have that 
$$\mathrm{ip}^{(2,X)@}_{\boldsymbol{\mathcal{A}}^{(2)}, s}(
\mathrm{CH}^{(2)}_{\boldsymbol{\mathcal{A}}^{(2)}, s} (
\mathfrak{P}^{(2), i, \bb{\mathfrak{P}^{(2)}}-1}
)
)
\mbox{ and }
\mathrm{ip}^{(2,X)@}_{\boldsymbol{\mathcal{A}}^{(2)}, s}(
\mathrm{CH}^{(2)}_{\boldsymbol{\mathcal{A}}^{(2)}, s}(
\mathfrak{P}^{(2),0,i-1}
)
)
$$
are paths in $\llbracket \mathfrak{P}^{(2),i, \bb{\mathfrak{P}^{(2)}}-1}\rrbracket_{s}$ and $\llbracket\mathfrak{P}^{(2),0,i-1}\rrbracket_{s}$, respectively. Hence the interpretation of the $1$-composition operation symbol $\circ^{1}_{s}$ in $\mathbf{F}_{\Sigma^{\boldsymbol{\mathcal{A}}^{(2)}}}(\mathbf{Pth}_{\boldsymbol{\mathcal{A}}^{(2)}})$ becomes that of $\mathbf{Pht}_{\boldsymbol{\mathcal{A}}^{(2)}}$.
The fourth equality holds since, according to Definition~\ref{DDIp}, $\mathrm{ip}^{(2,X)@}_{\boldsymbol{\mathcal{A}}^{(2)}}$ is a many-sorted $\Sigma^{\boldsymbol{\mathcal{A}}^{(2)}}$-homomorphism;
finally, the last equality recovers the definition of the second-order Curry-Howard mapping at a second-order path of length strictly greater than one with a second-order echelon at position $i$, introduced in Definition~\ref{DDCH}.

If~(2), i.e., if $\mathfrak{P}^{(2)}$ is an echelonless second-order path, it could be the case that~(2.1) $\mathfrak{P}^{(2)}$ is an echelonless second-order path that is not head-constant, or~(2.2) $\mathfrak{P}^{(2)}$ is a head-constant echelonless second-order path that is not coherent or~(2.3) $\mathfrak{P}^{(2)}$ is a coherent head-constant echelonless second-order path.

If~(2.1), i.e., if $\mathfrak{P}^{(2)}$ is an echelonless second-order path that is not head-constant, then, we let $i \in \bb{\mathfrak{P}^{(2)}}$ be the greatest index for which $\mathfrak{P}^{(2),0,i}$ is a head-constant echelonless second-order path. The following chain of equalities holds
\begin{flushleft}
$
\mathrm{ech}^{(2,\mathcal{A}^{(2)})\flat}_{\boldsymbol{\mathcal{A}}^{(2)}, s}(\mathfrak{P}^{(2)})
$
\allowdisplaybreaks
\begin{align*}
&=
\mathrm{ech}^{(2,\mathcal{A}^{(2)})\flat}_{\boldsymbol{\mathcal{A}}^{(2)}, s}\left(
\mathfrak{P}^{(2),i+1, \bb{\mathfrak{P}^{(2)}}-1}
\right)
\circ^{1\mathbf{Pth}_{\boldsymbol{\mathcal{A}}^{(1)}}}_{s}
\mathrm{ech}^{(2,\mathcal{A}^{(2)})\flat}_{\boldsymbol{\mathcal{A}}^{(2)}, s}\left(
\mathfrak{P}^{(2),0,i}
\right)
\tag{1}
\\
&=
\resizebox{.89\textwidth}{!}{$
\mathrm{ip}^{(2,X)@}_{\boldsymbol{\mathcal{A}}^{(2)}, s}\left(
\mathrm{CH}^{(2)}_{\boldsymbol{\mathcal{A}}^{(2)}, s} \left(
\mathfrak{P}^{(2), i+1, \bb{\mathfrak{P}^{(2)}}-1}
\right)
\right)
\circ^{1\mathbf{Pth}_{\boldsymbol{\mathcal{A}}^{(1)}}}_{s}
\mathrm{ip}^{(2,X)@}_{\boldsymbol{\mathcal{A}}^{(2)}, s}\left(
\mathrm{CH}^{(2)}_{\boldsymbol{\mathcal{A}}^{(2)}, s} \left(
\mathfrak{P}^{(2),0,i}
\right)
\right)
$}
\tag{2}
\\
&=
\resizebox{.89\textwidth}{!}{$
\mathrm{ip}^{(2,X)@}_{\boldsymbol{\mathcal{A}}^{(2)}, s}\left(
\mathrm{CH}^{(2)}_{\boldsymbol{\mathcal{A}}^{(2)}, s} \left(
\mathfrak{P}^{(2), i+1, \bb{\mathfrak{P}^{(2)}}-1}
\right)
\right)
\circ^{1\mathbf{F}_{\Sigma^{\boldsymbol{\mathcal{A}}^{(2)}}}(\mathbf{Pth}_{\boldsymbol{\mathcal{A}}^{(2)}})}_{s}
\mathrm{ip}^{(2,X)@}_{\boldsymbol{\mathcal{A}}^{(2)}, s}\left(
\mathrm{CH}^{(2)}_{\boldsymbol{\mathcal{A}}^{(2)}, s} \left(
\mathfrak{P}^{(2),0,i}
\right)
\right)
$}
\tag{3}
\\
&=
\mathrm{ip}^{(2,X)@}_{\boldsymbol{\mathcal{A}}^{(2)}, s}\left(
\mathrm{CH}^{(2)}_{\boldsymbol{\mathcal{A}}^{(2)}, s} \left(
\mathfrak{P}^{(2),i+1, \bb{\mathfrak{P}^{(2)}}-1}
\right)
\circ^{1\mathbf{T}_{\Sigma^{\boldsymbol{\mathcal{A}}^{(2)}}}(X)}_{s}
\mathrm{CH}^{(2)}_{\boldsymbol{\mathcal{A}}^{(2)}, s} \left(
\mathfrak{P}^{(2),0,i}
\right)
\right)
\tag{4}
\\
&=
\mathrm{ip}^{(2,X)@}_{\boldsymbol{\mathcal{A}}^{(2)}, s}\left(
\mathrm{CH}^{(2)}_{\boldsymbol{\mathcal{A}}^{(2)}, s} \left(
\mathfrak{P}^{(2)}
\right)
\right).
\tag{5}
\end{align*}
\end{flushleft}

In the just stated chain of equalities, the first equality unravels the definition of the second-order path extension mapping at an echelonless second-order path that is not head-constant, introduced in Proposition~\ref{PDPthExt};
the second equality follows taking into account that $(\mathfrak{P}^{((2),i,\bb{\mathfrak{P}^{(2)}}-1}, s)$ and $(\mathfrak{P}^{(2),0,i-1}, s)$ are strictly smaller that $(\mathfrak{P}^{(2)}, s)$ with respect to $\prec_{\mathbf{Pth}_{\boldsymbol{\mathcal{A}}^{(2)}}}$;
the third equality holds because, according to Proposition~\ref{PDIpDCH}, we have that 
$$\mathrm{ip}^{(2,X)@}_{\boldsymbol{\mathcal{A}}^{(2)}, s}(
\mathrm{CH}^{(2)}_{\boldsymbol{\mathcal{A}}^{(2)}, s} (
\mathfrak{P}^{(2), i+1, \bb{\mathfrak{P}^{(2)}}-1}
)
)
\mbox{ and }
\mathrm{ip}^{(2,X)@}_{\boldsymbol{\mathcal{A}}^{(2)}, s}(
\mathrm{CH}^{(2)}_{\boldsymbol{\mathcal{A}}^{(2)}, s}(
\mathfrak{P}^{(2),0,i}
)
)
$$
are paths in $\llbracket \mathfrak{P}^{(2),i+1, \bb{\mathfrak{P}^{(2)}}-1}\rrbracket_{s}$ and $\llbracket\mathfrak{P}^{(2),0,i}\rrbracket_{s}$, respectively. Hence the interpretation of the $1$-composition operation symbol $\circ^{1}_{s}$ in $\mathbf{F}_{\Sigma^{\boldsymbol{\mathcal{A}}^{(2)}}}(\mathbf{Pth}_{\boldsymbol{\mathcal{A}}^{(2)}})$ becomes that of $\mathbf{Pht}_{\boldsymbol{\mathcal{A}}^{(2)}}$.
The fourth equality holds since, according to Definition~\ref{DDIp}, $\mathrm{ip}^{(2,X)@}_{\boldsymbol{\mathcal{A}}^{(2)}}$ is a many-sorted $\Sigma^{\boldsymbol{\mathcal{A}}^{(2)}}$-homomorphism;
finally, the last equality recovers the definition of the second-order Curry-Howard mapping at an echelonless second-order path that is not head-constant, introduced in Definition~\ref{DDCH}.

If~(2.2), i.e., if $\mathfrak{P}^{(2)}$ is a head-constant echelonless second-order path that is not coherent, then, we let $i \in \bb{\mathfrak{P}^{(2)}}$ be the greatest index for which $\mathfrak{P}^{(2),0,i}$ is a head-constant echelonless coherent second-order path. The following chain of equalities holds
\begin{flushleft}
$
\mathrm{ech}^{(2,\mathcal{A}^{(2)})\flat}_{\boldsymbol{\mathcal{A}}^{(2)}, s}(\mathfrak{P}^{(2)})
$
\allowdisplaybreaks
\begin{align*}
&=
\mathrm{ech}^{(2,\mathcal{A}^{(2)})\flat}_{\boldsymbol{\mathcal{A}}^{(2)}, s}\left(
\mathfrak{P}^{(2),i+1, \bb{\mathfrak{P}^{(2)}}-1}
\right)
\circ^{1\mathbf{Pth}_{\boldsymbol{\mathcal{A}}^{(1)}}}_{s}
\mathrm{ech}^{(2,\mathcal{A}^{(2)})\flat}_{\boldsymbol{\mathcal{A}}^{(2)}, s}\left(
\mathfrak{P}^{(2),0,i}
\right)
\tag{1}
\\
&=
\resizebox{.89\textwidth}{!}{$
\mathrm{ip}^{(2,X)@}_{\boldsymbol{\mathcal{A}}^{(2)}, s}\left(
\mathrm{CH}^{(2)}_{\boldsymbol{\mathcal{A}}^{(2)}, s} \left(
\mathfrak{P}^{(2), i+1, \bb{\mathfrak{P}^{(2)}}-1}
\right)
\right)
\circ^{1\mathbf{Pth}_{\boldsymbol{\mathcal{A}}^{(1)}}}_{s}
\mathrm{ip}^{(2,X)@}_{\boldsymbol{\mathcal{A}}^{(2)}, s}\left(
\mathrm{CH}^{(2)}_{\boldsymbol{\mathcal{A}}^{(2)}, s} \left(
\mathfrak{P}^{(2),0,i}
\right)
\right)
$}
\tag{2}
\\
&=
\resizebox{.89\textwidth}{!}{$
\mathrm{ip}^{(2,X)@}_{\boldsymbol{\mathcal{A}}^{(2)}, s}\left(
\mathrm{CH}^{(2)}_{\boldsymbol{\mathcal{A}}^{(2)}, s} \left(
\mathfrak{P}^{(2), i+1, \bb{\mathfrak{P}^{(2)}}-1}
\right)
\right)
\circ^{1\mathbf{F}_{\Sigma^{\boldsymbol{\mathcal{A}}^{(2)}}}(\mathbf{Pth}_{\boldsymbol{\mathcal{A}}^{(2)}})}_{s}
\mathrm{ip}^{(2,X)@}_{\boldsymbol{\mathcal{A}}^{(2)}, s}\left(
\mathrm{CH}^{(2)}_{\boldsymbol{\mathcal{A}}^{(2)}, s} \left(
\mathfrak{P}^{(2),0,i}
\right)
\right)
$}
\tag{3}
\\
&=
\mathrm{ip}^{(2,X)@}_{\boldsymbol{\mathcal{A}}^{(2)}, s}\left(
\mathrm{CH}^{(2)}_{\boldsymbol{\mathcal{A}}^{(2)}, s} \left(
\mathfrak{P}^{(2),i+1, \bb{\mathfrak{P}^{(2)}}-1}
\right)
\circ^{1\mathbf{T}_{\Sigma^{\boldsymbol{\mathcal{A}}^{(2)}}}(X)}_{s}
\mathrm{CH}^{(2)}_{\boldsymbol{\mathcal{A}}^{(2)}, s} \left(
\mathfrak{P}^{(2),0,i}
\right)
\right)
\tag{4}
\\
&=
\mathrm{ip}^{(2,X)@}_{\boldsymbol{\mathcal{A}}^{(2)}, s}\left(
\mathrm{CH}^{(2)}_{\boldsymbol{\mathcal{A}}^{(2)}, s} \left(
\mathfrak{P}^{(2)}
\right)
\right).
\tag{5}
\end{align*}
\end{flushleft}

In the just stated chain of equalities, the first equality unravels the definition of the second-order path extension mapping at a head-constant echelonless second-order path that is not coherent, introduced in Proposition~\ref{PDPthExt};
the second equality follows taking into account that $(\mathfrak{P}^{((2),i,\bb{\mathfrak{P}^{(2)}}-1}, s)$ and $(\mathfrak{P}^{(2),0,i-1}, s)$ are strictly smaller that $(\mathfrak{P}^{(2)}, s)$ with respect to $\prec_{\mathbf{Pth}_{\boldsymbol{\mathcal{A}}^{(2)}}}$;
the third equality holds because, according to Proposition~\ref{PDIpDCH}, we have that 
$$\mathrm{ip}^{(2,X)@}_{\boldsymbol{\mathcal{A}}^{(2)}, s}(
\mathrm{CH}^{(2)}_{\boldsymbol{\mathcal{A}}^{(2)}, s} (
\mathfrak{P}^{(2), i+1, \bb{\mathfrak{P}^{(2)}}-1}
)
)
\mbox{ and }
\mathrm{ip}^{(2,X)@}_{\boldsymbol{\mathcal{A}}^{(2)}, s}(
\mathrm{CH}^{(2)}_{\boldsymbol{\mathcal{A}}^{(2)}, s}(
\mathfrak{P}^{(2),0,i}
)
)
$$
are paths in $\llbracket \mathfrak{P}^{(2),i+1, \bb{\mathfrak{P}^{(2)}}-1}\rrbracket_{s}$ and $\llbracket\mathfrak{P}^{(2),0,i}\rrbracket_{s}$, respectively. Hence the interpretation of the $1$-composition operation symbol $\circ^{1}_{s}$ in $\mathbf{F}_{\Sigma^{\boldsymbol{\mathcal{A}}^{(2)}}}(\mathbf{Pth}_{\boldsymbol{\mathcal{A}}^{(2)}})$ becomes that of $\mathbf{Pht}_{\boldsymbol{\mathcal{A}}^{(2)}}$.
The fourth equality holds since, according to Definition~\ref{DDIp}, $\mathrm{ip}^{(2,X)@}_{\boldsymbol{\mathcal{A}}^{(2)}}$ is a many-sorted $\Sigma^{\boldsymbol{\mathcal{A}}^{(2)}}$-homomorphism;
finally, the last equality recovers the definition of the second-order Curry-Howard mapping at a head-constant echelonless second-order path that is not coherent, introduced in Definition~\ref{DDCH}.

If~(2.3), i.e., if $\mathfrak{P}$ is a coherent head-constant echelonless path then, according to Definition~\ref{DDHeadCt}, there exsits a unique word $s \in S^{\star}-\{\lambda\}$ and a unique operation symbol $\tau$ in $\Sigma^{\boldsymbol{\mathcal{A}}^{(1)}}_{\mathbf{s}, s}$ associated to $\mathfrak{P}^{(2)}$. Let us recall from Defintion~\ref{DDPth} that the operations of $0$-source and $0$-target are forbidden. Let $(\mathfrak{P}^{(2)}_{j})_{j \in \bb{\mathbf{s}}}$ be the family of second-order paths we can extract from $\mathfrak{P}^{(2)}$ in virtue of Lemma~\ref{LDPthExtract}. Then, the following chain of equalities holds
\allowdisplaybreaks
\begin{align*}
\mathrm{ech}^{(2,\mathcal{A}^{(2)})\flat}_{\boldsymbol{\mathcal{A}}^{(2)}, s}(\mathfrak{P}^{(2)})
&=
\tau^{\mathbf{Pth}_{\boldsymbol{\mathcal{A}}^{(1)}}^{\mathrm{id}^{\boldsymbol{\mathcal{A}}^{(1)}}(1,2)}}\left(
\left(
\mathrm{ech}^{(2,\mathcal{A}^{(2)})\flat}_{\boldsymbol{\mathcal{A}}^{(2)}, s_{j}}(\mathfrak{P}^{(2)}_{j})
\right)_{j \in \bb{\mathbf{s}}}
\right)
\tag{1}
\\
&=
\tau^{\mathbf{Pth}_{\boldsymbol{\mathcal{A}}^{(2)}}^{\mathrm{id}^{\boldsymbol{\mathcal{A}}^{(1)}}(1,2)}}\left(
\left(
\mathrm{ip}^{(2,X)@}_{\boldsymbol{\mathcal{A}}^{(2)}, s}\left(
\mathrm{CH}^{(2)}_{\boldsymbol{\mathcal{A}}^{(2)}, s} \left(
\mathfrak{P}^{(2)}_{j}
\right)
\right)
\right)_{j \in \bb{\mathbf{s}}}
\right)
\tag{2}
\\
&=
\tau^{\mathbf{Pth}_{\boldsymbol{\mathcal{A}}^{(2)}}^{(1,2)}}\left(
\left(
\mathrm{ip}^{(2,X)@}_{\boldsymbol{\mathcal{A}}^{(2)}, s}\left(
\mathrm{CH}^{(2)}_{\boldsymbol{\mathcal{A}}^{(2)}, s} \left(
\mathfrak{P}^{(2)}_{j}
\right)
\right)
\right)_{j \in \bb{\mathbf{s}}}
\right)
\tag{3}
\\
&=
\tau^{\mathbf{F}_{\Sigma^{\boldsymbol{\mathcal{A}}^{(2)}}}(\mathbf{Pth}_{\boldsymbol{\mathcal{A}}^{(2)}})}\left(
\left(
\mathrm{ip}^{(2,X)@}_{\boldsymbol{\mathcal{A}}^{(2)}, s}\left(
\mathrm{CH}^{(2)}_{\boldsymbol{\mathcal{A}}^{(2)}, s} \left(
\mathfrak{P}^{(2)}_{j}
\right)
\right)
\right)_{j \in \bb{\mathbf{s}}}
\right)
\tag{4}
\\
&=
\mathrm{ip}^{(2,X)@}_{\boldsymbol{\mathcal{A}}^{(2)}, s}\left(
\tau^{\mathbf{T}_{\Sigma^{\boldsymbol{\mathcal{A}}^{(2)}}}(X)}\left(
\left(
\mathrm{CH}^{(2)}_{\boldsymbol{\mathcal{A}}^{(2)}, s} \left(
\mathfrak{P}^{(2)}_{j}
\right)
\right)
\right)_{j \in \bb{\mathbf{s}}}
\right)
\tag{5}
\\
&=
\mathrm{ip}^{(2,X)@}_{\boldsymbol{\mathcal{A}}^{(2)}, s}\left(
\mathrm{CH}^{(2)}_{\boldsymbol{\mathcal{A}}^{(2)}, s} \left(
\mathfrak{P}^{(2)}
\right)
\right).
\tag{6}
\end{align*}

In the just stated chain of equalities, the first equality unravels the definition of the path extension mapping at a coherent head-constant echelonless path, introduced in Proposition~\ref{PDPthExt};
the second equality follows taking into account that, for every $j\in\bb{\mathbf{s}}$, $(\mathfrak{P}^{(2)}_{j}, s_{j})$ is strictly smaller than $(\mathfrak{P}^{(2)}, s)$ with respect to $\prec_{\mathbf{Pth}_{\boldsymbol{\mathcal{A}}^{(2)}}}$;
the third equality follows from the fact that, according to Proposition~\ref{PAlgFun}, $\mathrm{Alg}_{\mathfrak{d}}$ is a functor, thus $\mathbf{Pth}_{\boldsymbol{\mathcal{A}}^{(2)}}^{\mathrm{id}^{\boldsymbol{\mathcal{A}}^{(2)}}(1,2)} 
= 
(\mathrm{id}^{\boldsymbol{\Sigma}})^{*}_{\mathfrak{d}}(\mathbf{Pth}_{\boldsymbol{\mathcal{A}}^{(2)}}^{(1,2)})
=
\mathbf{Pth}_{\boldsymbol{\mathcal{A}}^{(2)}}^{(1,2)}$;
the fourth equality holds because, according to Proposition~\ref{PDIpDCH}, we have that, for every $j \in \bb{\mathbf{s}}$,
$\mathrm{ip}^{(2,X)@}_{\boldsymbol{\mathcal{A}}^{(2)}, s}(
\mathrm{CH}^{(2)}_{\boldsymbol{\mathcal{A}}^{(2)}, s} (
\mathfrak{P}^{(2)}_{j}
)
)$
is a path in $\llbracket\mathfrak{P}^{(2)}_{j}\rrbracket_{s_{j}}$. Hence the interpretation of the operation symbol $\tau$ in $\mathbf{F}_{\Sigma^{\boldsymbol{\mathcal{A}}^{(2)}}}(\mathbf{Pth}_{\boldsymbol{\mathcal{A}}^{(2)}})$ becomes that of $\mathbf{Pht}_{\boldsymbol{\mathcal{A}}^{(2)}}$
the fifth equality holds since, according to Definition~\ref{DDIp}, $\mathrm{ip}^{(2,X)@}_{\boldsymbol{\mathcal{A}}^{(2)}}$ is a many-sorted $\Sigma^{\boldsymbol{\mathcal{A}}^{(2)}}$-homomorphism;
finally, the last equality recovers the definition of the second-order Curry-Howard mapping at a coherent head-constant echelonless path, introduced in Definition~\ref{DDCH}.

This proves the inductive step.

This completes the proof.
\end{proof}

\begin{proposition}
\label{PDQPthExtEch}
Let $\boldsymbol{\mathcal{A}}^{(2)}$ be a second-order rewriting system and let $\mathrm{id}^{\boldsymbol{\mathcal{A}}^{(2)}}$ be its second-order identity morphism. Thus, $\mathrm{ech}^{(2,\mathcal{A}^{(2)})@}_{\boldsymbol{\mathcal{A}}^{(2)}} = \mathrm{id}^{\llbracket\mathrm{Pth}_{\boldsymbol{\mathcal{A}}^{(2)}}\rrbracket}$.
\end{proposition}

\begin{proof}
Note that the following chain of equalities holds
\begin{align*}
\mathrm{id}^{\llbracket\mathbf{Pth}_{\boldsymbol{\mathcal{A}}^{(2)}}\rrbracket}
\circ
\mathrm{pr}^{\llbracket\cdot\rrbracket}_{\boldsymbol{\mathcal{A}}^{(2)}}
&=
\mathrm{ip}^{(\llbracket2\rrbracket,X)@}_{\boldsymbol{\mathcal{A}}^{(2)}}
\circ
\mathrm{CH}^{\llbracket2\rrbracket}_{\boldsymbol{\mathcal{A}}^{(2)}}
\circ
\mathrm{pr}^{\llbracket\cdot\rrbracket}_{\boldsymbol{\mathcal{A}}^{(2)}}
\tag{1}
\\
&=
\mathrm{ip}^{(\llbracket2\rrbracket,X)@}_{\boldsymbol{\mathcal{A}}^{(2)}}
\circ
\mathrm{pr}^{\Theta^{\llbracket2\rrbracket}_{\boldsymbol{\mathcal{A}}^{(2)}}}
\circ
\mathrm{CH}^{(2)\mathrm{m}}_{\boldsymbol{\mathcal{A}}^{(2)}}
\circ
\mathrm{pr}^{\llbracket\cdot\rrbracket}_{\boldsymbol{\mathcal{A}}^{(2)}}
\tag{2}
\\
&=
\mathrm{ip}^{(\llbracket2\rrbracket,X)@}_{\boldsymbol{\mathcal{A}}^{(2)}}
\circ
\mathrm{pr}^{\Theta^{\llbracket2\rrbracket}_{\boldsymbol{\mathcal{A}}^{(2)}}}
\circ
\mathrm{CH}^{(2)}_{\boldsymbol{\mathcal{A}}^{(2)}}
\tag{3}
\\
&=
\mathrm{pr}^{\llbracket\cdot\rrbracket}_{\boldsymbol{\mathcal{A}}^{(2)}}
\circ
\mathrm{ip}^{(2,X)@}_{\boldsymbol{\mathcal{A}}^{(2)}}
\circ
\mathrm{CH}^{(2)}_{\boldsymbol{\mathcal{A}}^{(2)}}
\tag{4}
\end{align*}

The first equality follows from Theorem~\ref{TDIso};
the second equality follows from the definition of the \(S\)-sorted mapping \(\mathrm{CH}^{\llbracket2\rrbracket}_{\boldsymbol{\mathcal{A}}^{(2)}}\), introduced in Definition~\ref{DDPTQDCH};
the third equality follows from the definition of the monomorphic Curry-Howard mapping \(\mathrm{CH}^{(2)\mathrm{m}}_{\boldsymbol{\mathcal{A}}^{(2)}}\), introduced in Definition~\ref{DDCHQuot};
finally, the last equality follows from the definition of the \(S\)-sorted mapping \(\mathrm{ip}^{(\llbracket2\rrbracket,X)@}_{\boldsymbol{\mathcal{A}}^{(2)}}\), introduced in Definition~\ref{DDPTQIp}.

Thus, by uniqueness of the second-order quotient path extension mapping introduced in Definition~\ref{DDQPthExt}, \(\mathrm{ech}^{(2,\boldsymbol{\mathcal{A}}^{(2)})@}\) is equal to \(\mathrm{id}^{\llbracket\mathrm{Pth}_{\boldsymbol{\mathcal{A}}^{(2)}}\rrbracket}\).
\end{proof}

Now we define the composition between second-order morphisms and show that it is not associative.

\begin{definition}
\label{DCompRws2}
Let $\mathbf{f}^{(2)}=(\varphi, c, (f^{(i)})_{i\in 3})$ be a second-order morphism from $\boldsymbol{\mathcal{A}}^{(2)}$ to $\boldsymbol{\mathcal{B}}^{(2)}$ and $\mathbf{g}^{(2)}=(\psi, d, (g^{(i)})_{i\in 3})$ a second-order morphism from $\boldsymbol{\mathcal{B}}^{(2)}$ to $\boldsymbol{\mathcal{C}}^{(2)}$. Its second-order \emph{composition} morphism, from $\boldsymbol{\mathcal{A}}^{(2)}$ to $\boldsymbol{\mathcal{C}}^{(2)}$, is given by
$$
\mathbf{g}^{(2)} \circ \mathbf{f}^{(2)}
=
\left(
\mathbf{g}^{(1)} \circ \mathbf{f}^{(1)},
g^{(2)\flat}_{\varphi} \circ f^{(2)}
\right).
$$
Note that $\mathbf{g}^{(1)} \circ \mathbf{f}^{(1)}$ is the first-order composition morphism, introduced in Definition~\ref{DCompRws1}, and $g^{(2)\flat}_{\varphi} \circ f^{(2)}$ is obtained from the following diagrams,
$$
\xymatrix{
\mathcal{B}^{(2)}
  \ar[r]^-{\mathrm{ech}^{(2,\mathcal{B}^{(2)})}_{\boldsymbol{\mathcal{B}}^{(2)}}}
  \ar[rd]_-{g^{(2)}}
&
\mathrm{Pth}_{\boldsymbol{\mathcal{B}}^{(2)}}
  \ar[d]^-{{g}^{(2)\flat}}
&&
\mathrm{Pth}_{\boldsymbol{\mathcal{B}}^{(2)}, \varphi}
  \ar[d]_-{g^{(2)\flat}_{\varphi}} 	
&
\mathcal{A^{(2)}}
  \ar[l]_-{f^{(2)}}
\\
& 
\mathrm{Pth}_{\boldsymbol{\mathcal{C}}^{(2)}, \psi}
&&
(\mathrm{Pth}_{\boldsymbol{\mathcal{C}}^{(2)}, \psi})_{\varphi}
\\
}
$$
being $g^{(2)\flat}$ the second-order path extension mapping of $g^{(2)}$ introduced in Proposition~\ref{PDPthExt}.
\end{definition}


\begin{remark}
The composition of second-order morphisms is not associative.
\end{remark}

We now define the notion of second-order equivalence of second-order morphisms. Thus, the composition of second-order morphisms will be associative up to second-order equivalence. That is, in order to define the category $\mathsf{Rws}_{\mathfrak{d}}^{\llbracket 2 \rrbracket}$, we will consider second-order morphisms up to second-order equivalence, i.e., equivalence classes of second-order morphisms.

\begin{definition}
\label{DDMorEqv}
Let $\mathbf{f}^{(2)}=(\varphi, c, (f^{(i)})_{i\in 3})$ and $\mathbf{g}^{(2)}=(\psi, d, (g^{(i)})_{i\in 3})$ be two second-order morphisms from $\boldsymbol{\mathcal{A}}^{(2)}$ to $\boldsymbol{\mathcal{B}}^{(2)}$. We will say that $\mathbf{f}^{(2)}$ and $\mathbf{g}^{(2)}$ are \emph{second-order equivalent}, written $\mathbf{f}^{(2)} \cong^{(2)} \mathbf{g}^{(2)}$, if 
\begin{enumerate}
\item
$(\varphi, c, f^{(0)}) = (\psi, d, g^{(0)})$, i.e., $\varphi = \psi$, $c=d$ and $f^{(0)}=g^{(0)}$; and
\item
$\mathrm{pr}^{\mathrm{Ker}(\mathrm{CH}^{(1)})}_{\boldsymbol{\mathcal{B}}^{(1)}, \varphi} \circ f^{(1)\flat} = \mathrm{pr}^{\mathrm{Ker}(\mathrm{CH}^{(1)})}_{\boldsymbol{\mathcal{B}}^{(1)}, \varphi} \circ g^{(1)\flat}$. That is, for every sort $s$ in $S$ and every path $\mathfrak{P}$ in $\mathrm{Pth}_{\boldsymbol{\mathcal{A}}, s}$,
$$
\left[
f_{s}^{(1)\flat}\left(
\mathfrak{P}
\right)
\right]_{\varphi(s)}
=
\left[
g_{s}^{(1)\flat}\left(
\mathfrak{P}
\right)
\right]_{\varphi(s)};
$$
and
\item
$\mathrm{pr}^{\llbracket\cdot\rrbracket}_{\boldsymbol{\mathcal{B}}^{(2)}, \varphi} \circ f^{(2)\flat} = \mathrm{pr}^{\llbracket\cdot\rrbracket}_{\boldsymbol{\mathcal{B}}^{(2)}, \varphi} \circ g^{(2)\flat}$. That is, for every sort $s$ in $S$ and every second-order path $\mathfrak{P}^{(2)}$ in $\mathrm{Pth}_{\boldsymbol{\mathcal{A}}^{(2)}, s}$,
$$
\left\llbracket
f_{s}^{(2)\flat}\left(
\mathfrak{P}^{(2)}
\right)
\right\rrbracket_{\varphi(s)}
=
\left\llbracket
g_{s}^{(2)\flat}\left(
\mathfrak{P}^{(2)}
\right)
\right\rrbracket_{\varphi(s)}.
$$
\end{enumerate}
Note that, items~(1) and (2) above imply that  $\mathbf{f}^{(1)}$ and $\mathbf{g}^{(1)}$ are first-order equivalent. Moreover, note that $\cong^{(2)}$ is an equivalence relation. Therefore, to simplify notation, we will denote by $\llbracket\mathbf{f}^{(2)}\rrbracket$ the equivalence class $[\mathbf{f}^{(2)}]_{\cong^{(2)}}$.
\end{definition}

In the next proposition we state the relation between the path-extension mapping of a composition of second-order morphisms and the path-extension mappings of each of the composites. Let us say that, contrary to what happens with the homomorphic extension at level 0, it is not necessarily true that $(g^{(2)\flat}_{\varphi} \circ  f^{(2)} )^{\flat}$ is equal to $g^{(2)\flat}_{\varphi} \circ f^{(2)\flat}$.

\begin{proposition}
\label{PDPthExtComp}
Let $\mathbf{f}^{(2)}=(\varphi, c, (f^{(i)})_{i\in 3})$ be a second-order morphism from $\boldsymbol{\mathcal{A}}^{(2)}$ to $\boldsymbol{\mathcal{B}}^{(2)}$ and $\mathbf{g}^{(2)}=(\psi, d, (g^{(i)})_{i\in 3})$ a second-order morphism from $\boldsymbol{\mathcal{B}}^{(2)}$ to $\boldsymbol{\mathcal{C}}^{(2)}$. Then
$$
\mathrm{pr}^{\llbracket\cdot\rrbracket}_{\boldsymbol{\mathcal{C}}^{(2)}, \psi \circ \varphi}
\circ
\left(
g^{(2)\flat}_{\varphi} \circ 
f^{(2)}
\right)^{\flat}
=
\mathrm{pr}^{\llbracket\cdot\rrbracket}_{\boldsymbol{\mathcal{C}}^{(2)}, \psi \circ \varphi}
\circ
g^{(2)\flat}_{\varphi} \circ
f^{(2)\flat}
$$
The reader is advised to consult the diagram of Figure~\ref{FDPthExtComp}.

\begin{figure}
$$
\xymatrix{
\mathcal{A}^{(2)}
	\ar[r]^-{\mathrm{ech}^{(2,\boldsymbol{\mathcal{A}}^{(2)})}}
	\ar[rdd]_-{g^{(2)\flat}_{\varphi} \circ f^{(2)}}
	&
\mathrm{Pth}_{\boldsymbol{\mathcal{A}}^{(2)}}
	\ar[dd]^-{\left(g^{(2)\flat}_{\varphi} \circ f^{(2)}\right)^{\flat}}
	&&
\mathrm{Pth}_{\boldsymbol{\mathcal{A}}^{(2)}}
	\ar[d]_-{f^{(2)\flat}}
	&
\mathcal{A}^{(2)}
	\ar[l]_-{\mathrm{ech}^{(2,\boldsymbol{\mathcal{A}}^{(2)})}}
	\ar[ld]^-{f^{(2)}}
\\
	&&&
\mathrm{Pth}_{\boldsymbol{\mathcal{B}}^{(2)}, \varphi}
	\ar[d]_-{g^{(2)\flat}_{\varphi}}
\\
	&
\mathrm{Pth}_{\boldsymbol{\mathcal{C}}^{(2)}, \psi\circ\varphi}
	\ar[d]^-{\mathrm{pr}^{\llbracket\cdot\rrbracket}_{\boldsymbol{\mathcal{C}}^{(2)}, \psi\circ\varphi}}
	&&
\mathrm{Pth}_{\boldsymbol{\mathcal{C}}^{(2)}, \psi\circ\varphi}
	\ar[d]_-{\mathrm{pr}^{\llbracket\cdot\rrbracket}_{\boldsymbol{\mathcal{C}}^{(2)}, \psi\circ\varphi}}
\\
	&
\llbracket\mathrm{Pth}_{\boldsymbol{\mathcal{C}}^{(2)}}\rrbracket_{\psi\circ\varphi}
	&&
\llbracket\mathrm{Pth}_{\boldsymbol{\mathcal{C}}^{(2)}}\rrbracket_{\psi\circ\varphi}
}
$$
\caption{Relation between the second-order path extension mapping of a composition and the corresponding second-order path extension mappings.}
\label{FDPthExtComp}
\end{figure}
\end{proposition}

\begin{proof}
Let $s$ in $S$ be a sort and $\mathfrak{P}^{(2)}$ a second-order path in $\mathrm{Pth}_{\boldsymbol{\mathcal{A}}^{(2)}, s}$. We will prove that
$$
\left\llbracket
\left(g^{(2)\flat}_{\varphi} \circ f^{(2)}\right)^{\flat}_{s}\left(
\mathfrak{P}^{(2)}
\right)
\right\rrbracket_{\psi(\varphi(s))}
=
\left\llbracket
g^{(2)\flat}_{\varphi(s)} \circ f^{(2)\flat}_{s}\left(
\mathfrak{P}^{(2)}
\right)
\right\rrbracket_{\psi(\varphi(s))}
$$
by Artinian induction on $(\coprod \mathrm{Pth}_{\boldsymbol{\mathcal{A}}^{(2)}}, \leq_{\mathbf{Pth}_{\boldsymbol{\mathcal{A}}^{(2)}}})$.

{\sffamily Base step of the Artinian induction.}

Let $(\mathfrak{P}^{(2)}, s)$ be a minimal element of $(\coprod \mathrm{Pth}_{\boldsymbol{\mathcal{A}}^{(2)}}, \leq_{\mathbf{Pth}_{\boldsymbol{\mathcal{A}}^{(2)}}})$. Then, by Proposition~\ref{PDMinimal}, the path $\mathfrak{P}^{(2)}$ is either~(1) an $(2,[1])$-identity second-order path or~(2) a second-order echelon.

If~(1), i.e., if $\mathfrak{P}^{(2)}$ is a $(2,[1])$-identity second-order path, then $\mathfrak{P}^{(2)}=\mathrm{ip}_{\boldsymbol{\mathcal{A}}^{(2)},s}^{(2,[1])\sharp}([P]_{s})$ for some path term class $[P]_{s} \in [\mathrm{PT}_{\boldsymbol{\mathcal{A}}^{(1)}}]_{s}$. Then the following chain of equalities holds
\allowdisplaybreaks
\begin{flalign*}
\left(
g^{(2)\flat}_{\varphi}
\circ
f^{(2)}
\right)^{\flat}_{s}\left(
\mathfrak{P}^{(2)}
\right)
&=
\left(
g^{(2)\flat}_{\varphi}
\circ
f^{(2)}
\right)^{\flat}_{s}\left(
\mathrm{ip}_{\boldsymbol{\mathcal{A}}^{(2)},s}^{(2,[1])\sharp}\left(
\left[P\right]_{s}
\right)
\right)
\tag{1}
\\
&=
\mathrm{ip}_{\boldsymbol{\mathcal{C}}^{(2)}, \psi(\varphi(s))}^{(2,[1])\sharp} \left(
\left(
g^{(1)\flat}_{\varphi} \circ f^{(1)}
\right)_{s}^{\mathsf{q}} \left(
\left[P\right]_{s}
\right)
\right)
\tag{2}
\\
&=
\mathrm{ip}_{\boldsymbol{\mathcal{C}}^{(2)}, \psi(\varphi(s))}^{(2,[1])\sharp} \left(
g^{(1)@}_{\varphi(s)} \left(
f_{s}^{(1)@} \left(
\left[P\right]_{s}
\right) 
\right)
\right)
\tag{3}
\\
&=
g_{\varphi(s)}^{(2)\flat} \left(
\mathrm{ip}_{\boldsymbol{\mathcal{B}}^{(2)}, \varphi(s)}^{(2,[1])\sharp} \left(
f_{s}^{(1)@} \left(
\left[P\right]_{s}
\right) 
\right)
\right)
\tag{4}
\\
&=
g^{(2)\flat}_{\varphi(s)} \left(
f_{s}^{(2)\flat} \left(
\mathrm{ip}_{\boldsymbol{\mathcal{A}}^{(2)}, s}^{(2,[1])\sharp} \left(
\left[P\right]_{s}
\right) 
\right)
\right)
\tag{5}
\\
&=
g^{(2)\flat}_{\varphi(s)} \left(
f_{s}^{(2)\flat} \left(
\mathfrak{P}^{(2)}
\right)
\right).
\tag{6}
\end{flalign*}

The first equality unravels the definition of the second-order path $\mathfrak{P}^{(2)}$;
the second equality follows from Proposition~\ref{PDPthExt};
the third equality follows from Corollary~\ref{CQPthExtComp} and unraveling the $s$-th component of the composition of $S$-sorted mappings;
the fourth and fifth equalities follow from item~(2) of the definition of the second-order path extension mapping introduced in Proposition~\ref{PDPthExt};
finally, the last equality recovers the definition of the second-order path $\mathfrak{P}^{(2)}$.

Therefore, their $\llbracket \cdot \rrbracket$-classes coincide.

This completes case (1).

If~(2), i.e., if $\mathfrak{P}^{(2)}$ is a second-order echelon associated to a second-order rewrite rule $\mathfrak{p}^{(2)}=([M]_{s}, [N]_{s})$, that is, if $\mathfrak{P}^{(2)}$ has the form
$$
\xymatrix@C=105pt{
\mathfrak{P}^{(2)}: [M]_{s}
\ar@{=>}[r]^-{\text{\Small{$(\mathfrak{p}^{(2)},\mathrm{id}^{\mathrm{T}_{\Sigma^{\boldsymbol{\mathcal{A}}^{(1)}}}(X)_{s}})$}}}
&
[N]_{s}
},
$$

Then the following chain of equalities holds
\begin{flalign*}
\left(
g^{(2)\flat}_{\varphi}
\circ
f^{(2)}
\right)^{\flat}_{s}\left(
\mathfrak{P}^{(2)}
\right)
&=
\left(
g^{(2)\flat}_{\varphi}
\circ
f^{(2)}
\right)^{\flat}_{s}\left(
\mathrm{ech}_{\boldsymbol{\mathcal{A}}^{(2)}, s}^{(2,\mathcal{A}^{(2)})}\left(
\mathfrak{p}^{(2)}
\right)
\right)
\tag{1}
\\
&=
\left(
g^{(2)\flat}_{\varphi}
\circ
f^{(2)}
\right)_{s}\left(
\mathfrak{p}^{(2)}
\right)
\tag{2}
\\
&=
g^{(2)\flat}_{\varphi(s)} \left(
f_{s}^{(2)} \left(
\mathfrak{p}^{(2)}
\right)
\right)
\tag{3}
\\
&=
g^{(2)\flat}_{\varphi(s)} \left(
f_{s}^{(2)\flat} \left(
\mathrm{ech}_{\boldsymbol{\mathcal{A}}^{(2)}, s}^{(2,\mathcal{A}^{(2)})}\left(
\mathfrak{p}^{(2)}
\right)
\right)
\right)
\tag{4}
\\
&=
g^{(2)\flat}_{\varphi(s)} \left(
f_{s}^{(2)\flat} \left(
\mathfrak{P}^{(2)}
\right)
\right).
\tag{5}
\end{flalign*}

The first equality unravels the definition of the path $\mathfrak{P}^{(2)}$;
the second equality follows from Proposition~\ref{PDPthExt} considering the second-order morphism $\mathbf{g}^{(2)} \circ \mathbf{f}^{(2)}$;
the third equality unravels the definition of the $s$-th component of the composition of $S$-sorted mappings;
the fourth equality follows from Proposition~\ref{PDPthExt} considering the second-order morphism $\mathbf{f}^{(2)}$;
finally, the last equality recovers the definition of the path $\mathfrak{P}^{(2)}$.

Therefore, their $\llbracket \cdot \rrbracket$-classes coincide.

This completes case (2).

This concludes the base step of the Artinian induction.

{\sffamily Inductive step of the Artinian induction.}

Let $(\mathfrak{P}^{(2)},s)$ be a non-minimal element of $(\coprod\mathrm{Pth}_{\boldsymbol{\mathcal{A}}^{(2)}}, \leq_{\mathbf{Pth}_{\boldsymbol{\mathcal{A}}^{(2)}}})$. We can assume that $\mathfrak{P}^{(2)}$ is a not a $(2,[1])$-identity second-order path, since for those paths the desired equality has already been proven. Let us suppose that, for every sort $t\in S$ and every path $\mathfrak{Q}^{(2)}\in\mathrm{Pth}_{\boldsymbol{\mathcal{A}}^{(2)},t}$, if $(\mathfrak{Q}^{(2)},t)<_{\mathbf{Pth}_{\boldsymbol{\mathcal{A}}^{(2)}}}(\mathfrak{P}^{(2)},s)$, then the equality
$$
\left\llbracket
\left(g^{(2)\flat}_{\varphi} \circ f^{(2)}\right)^{\flat}_{t}\left(
\mathfrak{Q}^{(2)}
\right)
\right\rrbracket_{\psi(\varphi(t))}
=
\left\llbracket
g^{(2)\flat}_{\varphi(t)} \circ f^{(2)\flat}_{t}\left(
\mathfrak{Q}^{(2)}
\right)
\right\rrbracket_{\psi(\varphi(t))}
$$
holds.

By Lemma~\ref{LDOrdI}, we have that $\mathfrak{P}$ is either~(1) a second-order path of length strictly greater than one containing at least one second-order echelon or~(2) an
echelonless second-order path.

If~(1), i.e., if $\mathfrak{P}^{(2)}$ is a second-order path of length strictly greater than one containing at least one second-order echelon, then let $i\in \bb{\mathfrak{P}^{(2)}}$ be the first index for which the one-step subpath $\mathfrak{P}^{(2)i,i}$ of $\mathfrak{P}^{(2)}$ is a second-order echelon. We consider different cases for $i$ according to the cases presented in Definition~\ref{DDOrd}.

If~$i=0$, we have that the pairs $(\mathfrak{P}^{(2)0,0}, s)$ and $(\mathfrak{P}^{(2)1,\bb{\mathfrak{P}^{(2)}}-1}, s)$ $\prec_{\mathbf{Pth}_{\boldsymbol{\mathcal{A}}^{(2)}}}$-precede the pair $(\mathfrak{P}^{(2)}, s)$. Therefore, by induction, the equalities
$$
\left\llbracket
\left(g^{(2)\flat}_{\varphi} \circ f^{(2)}\right)^{\flat}_{s}\left(
\mathfrak{P}^{(2)0,0}
\right)
\right\rrbracket_{\psi(\varphi(s))}
=
\left\llbracket
g^{(2)\flat}_{\varphi(s)} \circ f^{(2)\flat}_{s}\left(
\mathfrak{P}^{(2)0,0}
\right)
\right\rrbracket_{\psi(\varphi(s))}
$$
$$
\left\llbracket
\left(g^{(2)\flat}_{\varphi} \circ f^{(2)}\right)^{\flat}_{s}\left(
\mathfrak{P}^{(2)1,\bb{\mathfrak{P}^{(2)}}-1}
\right)
\right\rrbracket_{\psi(\varphi(s))}
=
\left\llbracket
g^{(2)\flat}_{\varphi(s)} \circ f^{(2)\flat}_{s}\left(
\mathfrak{P}^{(2)1,\bb{\mathfrak{P}^{(2)}}-1}
\right)
\right\rrbracket_{\psi(\varphi(s))}
$$
hold. Then, the following chain of equalities holds
\begin{flushleft}
$
\left\llbracket
\left(
g^{(2)\flat}_{\varphi} \circ f^{(2)}
\right)_{s}^{\flat} \left(
\mathfrak{P}^{(2)}
\right)
\right\rrbracket_{\psi(\varphi(s))}
$
\allowdisplaybreaks
\begin{align*}
&=
\left\llbracket
\left(
g^{(2)\flat}_{\varphi} \circ f^{(2)}
\right)_{s}^{\flat} \left(
\mathfrak{P}^{(2)1, \bb{\mathfrak{P}^{(2)}}-1} \circ_{s}^{1\mathbf{Pth}_{\boldsymbol{\mathcal{A}}^{(2)}}} \mathfrak{P}^{(2)0,0}
\right)
\right\rrbracket_{\psi(\varphi(s))}
\tag{1}
\\
&=
\left\llbracket
\left(
g^{(2)\flat}_{\varphi} \circ f^{(2)}
\right)_{s}^{\flat} \left(
\mathfrak{P}^{(2)1, \bb{\mathfrak{P}^{(2)}}-1}
\right)
\right\rrbracket_{\psi(\varphi(s))}
\\
&\hspace{3.5cm}
\circ_{s}^{1\left\llbracket\mathbf{Pth}_{\boldsymbol{\mathcal{C}}^{(2)}}^{\mathbf{g}^{(2)}\circ\mathbf{f}^{(2)}}\right\rrbracket}
\left\llbracket
\left(
g^{(2)\flat}_{\varphi} \circ f^{(2)}
\right)_{s}^{\flat} \left(
\mathfrak{P}^{(2)0,0}
\right)
\right\rrbracket_{\psi(\varphi(s))}
\tag{2}
\\
&=
\left\llbracket
\left(g^{(2)\flat}_{\varphi} \circ f^{(2)\flat}\right)_{s} \left(
\mathfrak{P}^{(2)1, \bb{\mathfrak{P}^{(2)}}-1}
\right)
\right\rrbracket_{\psi(\varphi(s))}
\\
&\hspace{3.5cm}
\circ_{s}^{1\left\llbracket\mathbf{Pth}_{\boldsymbol{\mathcal{C}}^{(1)}}^{\mathbf{g}^{(2)}\circ\mathbf{f}^{(2)}}\right\rrbracket}
\left[
\left(g^{(2)\flat}_{\varphi} \circ f^{(2)\flat}\right)_{s} \left(
\mathfrak{P}^{(2)0,0}
\right)
\right\rrbracket_{\psi(\varphi(s))}
\tag{3}
\\
&=
\left\llbracket
g_{\varphi(s)}^{(2)\flat} \left(
f_{s}^{(2)\flat} \left(
\mathfrak{P}^{(2)1, \bb{\mathfrak{P}^{(2)}}-1}
\right)
\right)
\right\rrbracket_{\psi(\varphi(s))}
\\
&\hspace{3.5cm}
\circ_{\varphi(s)}^{1\llbracket\mathbf{Pth}_{\boldsymbol{\mathcal{C}}^{(2)}}^{\mathbf{g}^{(2)}}\rrbracket}
\left\llbracket
g_{\varphi(s)}^{(2)\flat} \left(
f_{s}^{(2)\flat} \left(
\mathfrak{P}^{(2)0,0}
\right)
\right)
\right\rrbracket_{\psi(\varphi(s))}
\tag{4}
\\
&=
\left\llbracket
g_{\varphi(s)}^{(2)\flat} \left(
f_{s}^{(2)\flat} \left(
\mathfrak{P}^{(2)1, \bb{\mathfrak{P}^{(2)}}-1}
\right)
\circ_{\varphi(s)}^{1\mathbf{Pth}_{\boldsymbol{\mathcal{B}}^{(2)}}}
f_{s}^{(2)\flat} \left(
\mathfrak{P}^{(2)0,0}
\right)
\right)
\right\rrbracket_{\psi(\varphi(s))}
\tag{5}
\\
&=
\left\llbracket
g_{\varphi(s)}^{(2)\flat} \left(
f_{s}^{(2)\flat} \left(
\mathfrak{P}^{(2)1, \bb{\mathfrak{P}^{(2)}}-1}
\circ_{s}^{1\mathbf{Pth}_{\boldsymbol{\mathcal{A}}^{(2)}}}
\mathfrak{P}^{(2)0,0}
\right)
\right)
\right\rrbracket_{\psi(\varphi(s))}
\tag{6}
\\
&=
\left\llbracket
g_{\varphi(s)}^{(2)\flat} \left(
f_{s}^{(2)\flat} \left(
\mathfrak{P}^{(2)}
\right)
\right)
\right\rrbracket_{\psi(\varphi(s))}.
\tag{7}
\end{align*}
\end{flushleft}

The first equality follows from the fact that $\mathfrak{P}^{(2)} = \mathfrak{P}^{(2)1, \bb{\mathfrak{P}^{(2)}}-1} \circ_{s}^{1\mathbf{Pth}_{\boldsymbol{\mathcal{A}}^{(2)}}} \mathfrak{P}^{(2)0,0}$;
the second equality follows from Proposition~\ref{PDHomPthExtKer};
the third equality follows by Artinian induction;
the fourth equality unravels the definition of the $s$-th component of the composition of $S$-sorted mappings and note that, for every sort $s$ in $S$,
$$
\circ_{s}^{1\llbracket\mathbf{Pth}_{\boldsymbol{\mathcal{C}}^{(2)}}^{\mathbf{g}^{(2)}\circ\mathbf{f}^{(2)}}\rrbracket}
=
\circ_{\psi(\varphi(s))}^{\llbracket\mathbf{Pth}_{\boldsymbol{\mathcal{C}}^{(2)}}\rrbracket}
=
\circ_{\varphi(s)}^{1\llbracket\mathbf{Pth}_{\boldsymbol{\mathcal{C}}^{(2)}}^{\mathbf{g}^{(2)}}\rrbracket};
$$
the fifth equality follows from Proposition~\ref{PDHomPthExtKer};
the sixth equality recovers the definition of $f^{(2)\flat}$;
finally, the last equality recovers the definition of $\mathfrak{P}^{(2)}$;

If~$i>0$, we have that the pairs $(\mathfrak{P}^{(2)0,i-1}, s)$ and $(\mathfrak{P}^{(2)i,\bb{\mathfrak{P}^{(2)}}-1}, s)$ $\prec_{\mathbf{Pth}_{\boldsymbol{\mathcal{A}}^{(2)}}}$-precede the pair $(\mathfrak{P}^{(2)}, s)$. Therefore, by induction, the equalities
$$
\left\llbracket
\left(g^{(2)\flat}_{\varphi} \circ f^{(2)}\right)^{\flat}_{s}\left(
\mathfrak{P}^{(2)0,i-1}
\right)
\right\rrbracket_{\psi(\varphi(s))}
=
\left\llbracket
g^{(2)\flat}_{\varphi(s)} \circ f^{(2)\flat}_{s}\left(
\mathfrak{P}^{(2)0,i-1}
\right)
\right\rrbracket_{\psi(\varphi(s))}
$$
$$
\left\llbracket
\left(g^{(2)\flat}_{\varphi} \circ f^{(2)}\right)^{\flat}_{s}\left(
\mathfrak{P}^{(2)i,\bb{\mathfrak{P}^{(2)}}-1}
\right)
\right\rrbracket_{\psi(\varphi(s))}
=
\left\llbracket
g^{(2)\flat}_{\varphi(s)} \circ f^{(2)\flat}_{s}\left(
\mathfrak{P}^{(2)i,\bb{\mathfrak{P}^{(2)}}-1}
\right)
\right\rrbracket_{\psi(\varphi(s))}
$$
hold. Then, the following chain of equalities holds
\begin{flushleft}
$
\left\llbracket
\left(
g^{(2)\flat}_{\varphi} \circ f^{(2)}
\right)_{s}^{\flat} \left(
\mathfrak{P}^{(2)}
\right)
\right\rrbracket_{\psi(\varphi(s))}
$
\allowdisplaybreaks
\begin{align*}
&=
\left\llbracket
\left(
g^{(2)\flat}_{\varphi} \circ f^{(2)}
\right)_{s}^{\flat} \left(
\mathfrak{P}^{(2)i, \bb{\mathfrak{P}^{(2)}}-1} \circ_{s}^{1\mathbf{Pth}_{\boldsymbol{\mathcal{A}}^{(2)}}} \mathfrak{P}^{(2)0,i-1}
\right)
\right\rrbracket_{\psi(\varphi(s))}
\tag{1}
\\
&=
\left\llbracket
\left(
g^{(2)\flat}_{\varphi} \circ f^{(2)}
\right)_{s}^{\flat} \left(
\mathfrak{P}^{(2)i, \bb{\mathfrak{P}^{(2)}}-1}
\right)
\right\rrbracket_{\psi(\varphi(s))}
\\
&\hspace{3.5cm}
\circ_{s}^{1\left\llbracket\mathbf{Pth}_{\boldsymbol{\mathcal{C}}^{(2)}}^{\mathbf{g}^{(2)}\circ\mathbf{f}^{(2)}}\right\rrbracket}
\left\llbracket
\left(
g^{(2)\flat}_{\varphi} \circ f^{(2)}
\right)_{s}^{\flat} \left(
\mathfrak{P}^{(2)0,i-1}
\right)
\right\rrbracket_{\psi(\varphi(s))}
\tag{2}
\\
&=
\left\llbracket
\left(g^{(2)\flat}_{\varphi} \circ f^{(2)\flat}\right)_{s} \left(
\mathfrak{P}^{(2)i, \bb{\mathfrak{P}^{(2)}}-1}
\right)
\right\rrbracket_{\psi(\varphi(s))}
\\
&\hspace{3.5cm}
\circ_{s}^{1\left\llbracket\mathbf{Pth}_{\boldsymbol{\mathcal{C}}^{(1)}}^{\mathbf{g}^{(2)}\circ\mathbf{f}^{(2)}}\right\rrbracket}
\left[
\left(g^{(2)\flat}_{\varphi} \circ f^{(2)\flat}\right)_{s} \left(
\mathfrak{P}^{(2)0,i-1}
\right)
\right\rrbracket_{\psi(\varphi(s))}
\tag{3}
\\
&=
\left\llbracket
g_{\varphi(s)}^{(2)\flat} \left(
f_{s}^{(2)\flat} \left(
\mathfrak{P}^{(2)i, \bb{\mathfrak{P}^{(2)}}-1}
\right)
\right)
\right\rrbracket_{\psi(\varphi(s))}
\\
&\hspace{3.5cm}
\circ_{\varphi(s)}^{1\llbracket\mathbf{Pth}_{\boldsymbol{\mathcal{C}}^{(2)}}^{\mathbf{g}^{(2)}}\rrbracket}
\left\llbracket
g_{\varphi(s)}^{(2)\flat} \left(
f_{s}^{(2)\flat} \left(
\mathfrak{P}^{(2)0,i-1}
\right)
\right)
\right\rrbracket_{\psi(\varphi(s))}
\tag{4}
\\
&=
\left\llbracket
g_{\varphi(s)}^{(2)\flat} \left(
f_{s}^{(2)\flat} \left(
\mathfrak{P}^{(2)i, \bb{\mathfrak{P}^{(2)}}-1}
\right)
\circ_{\varphi(s)}^{1\mathbf{Pth}_{\boldsymbol{\mathcal{B}}^{(2)}}}
f_{s}^{(2)\flat} \left(
\mathfrak{P}^{(2)0,i-1}
\right)
\right)
\right\rrbracket_{\psi(\varphi(s))}
\tag{5}
\\
&=
\left\llbracket
g_{\varphi(s)}^{(2)\flat} \left(
f_{s}^{(2)\flat} \left(
\mathfrak{P}^{(2)i, \bb{\mathfrak{P}^{(2)}}-1}
\circ_{s}^{1\mathbf{Pth}_{\boldsymbol{\mathcal{A}}^{(2)}}}
\mathfrak{P}^{(2)0,i-1}
\right)
\right)
\right\rrbracket_{\psi(\varphi(s))}
\tag{6}
\\
&=
\left\llbracket
g_{\varphi(s)}^{(2)\flat} \left(
f_{s}^{(2)\flat} \left(
\mathfrak{P}^{(2)}
\right)
\right)
\right\rrbracket_{\psi(\varphi(s))}.
\tag{7}
\end{align*}
\end{flushleft}

The first equality follows from the fact that $\mathfrak{P}^{(2)} = \mathfrak{P}^{(2)i, \bb{\mathfrak{P}^{(2)}}-1} \circ_{s}^{1\mathbf{Pth}_{\boldsymbol{\mathcal{A}}^{(2)}}} \mathfrak{P}^{(2)0,i-1}$;
the second equality follows from Proposition~\ref{PDHomPthExtKer};
the third equality follows by Artinian induction;
the fourth equality unravels the definition of the $s$-th component of the composition of $S$-sorted mappings and note that, for every sort $s$ in $S$,
$$
\circ_{s}^{1\llbracket\mathbf{Pth}_{\boldsymbol{\mathcal{C}}^{(2)}}^{\mathbf{g}^{(2)}\circ\mathbf{f}^{(2)}}\rrbracket}
=
\circ_{\psi(\varphi(s))}^{\llbracket\mathbf{Pth}_{\boldsymbol{\mathcal{C}}^{(2)}}\rrbracket}
=
\circ_{\varphi(s)}^{1\llbracket\mathbf{Pth}_{\boldsymbol{\mathcal{C}}^{(2)}}^{\mathbf{g}^{(2)}}\rrbracket};
$$
the fifth equality follows from Proposition~\ref{PDHomPthExtKer};
the sixth equality recovers the definition of $f^{(2)\flat}$;
finally, the last equality recovers the definition of $\mathfrak{P}^{(2)}$;

Case (1) follows.

If~(2), i.e., if $\mathfrak{P}^{(2)}$ is an echelonless second-order path in $\mathrm{Pth}_{\boldsymbol{\mathcal{A}}^{(2)},s}$. It could be the case that $(2.1)$ $\mathfrak{P}^{(2)}$ is not head-constant. Then let $i \in \bb{\mathfrak{P}^{(2)}}$ be the maximum index for which the subpath $\mathfrak{P}^{(2),0,i}$ of $\mathfrak{P}^{(2)}$ is a head-constant, echelonless second-order path. Note that the pairs $(\mathfrak{P}^{(2),0,i}, s)$ and $(\mathfrak{P}^{(2), i+1, \bb{\mathfrak{P}^{(2)}}-1}, s)$ $\prec_{\mathbf{Pth}_{\boldsymbol{\mathcal{A}}^{(2)}}}$-precede the pair $(\mathfrak{P}^{(2)},s)$. Therefore, by induction, the equalities
$$
\left\llbracket
\left(g^{(2)\flat}_{\varphi} \circ f^{(2)}\right)^{\flat}_{s}\left(
\mathfrak{P}^{(2)0,i}
\right)
\right\rrbracket_{\psi(\varphi(s))}
=
\left\llbracket
g^{(2)\flat}_{\varphi(s)} \circ f^{(2)\flat}_{s}\left(
\mathfrak{P}^{(2)0,i}
\right)
\right\rrbracket_{\psi(\varphi(s))}
$$
$$
\left\llbracket
\left(g^{(2)\flat}_{\varphi} \circ f^{(2)}\right)^{\flat}_{s}\left(
\mathfrak{P}^{(2)i+1,\bb{\mathfrak{P}^{(2)}}-1}
\right)
\right\rrbracket_{\psi(\varphi(s))}
=
\left\llbracket
g^{(2)\flat}_{\varphi(s)} \circ f^{(2)\flat}_{s}\left(
\mathfrak{P}^{(2)i+1,\bb{\mathfrak{P}^{(2)}}-1}
\right)
\right\rrbracket_{\psi(\varphi(s))}
$$
hold. Then, the following chain of equalities holds
\begin{flushleft}
$
\left\llbracket
\left(
g^{(2)\flat}_{\varphi} \circ f^{(2)}
\right)_{s}^{\flat} \left(
\mathfrak{P}^{(2)}
\right)
\right\rrbracket_{\psi(\varphi(s))}
$
\allowdisplaybreaks
\begin{align*}
&=
\left\llbracket
\left(
g^{(2)\flat}_{\varphi} \circ f^{(2)}
\right)_{s}^{\flat} \left(
\mathfrak{P}^{(2)i+1, \bb{\mathfrak{P}^{(2)}}-1} \circ_{s}^{1\mathbf{Pth}_{\boldsymbol{\mathcal{A}}^{(2)}}} \mathfrak{P}^{(2)0,i}
\right)
\right\rrbracket_{\psi(\varphi(s))}
\tag{1}
\\
&=
\left\llbracket
\left(
g^{(2)\flat}_{\varphi} \circ f^{(2)}
\right)_{s}^{\flat} \left(
\mathfrak{P}^{(2)i+1, \bb{\mathfrak{P}^{(2)}}-1}
\right)
\right\rrbracket_{\psi(\varphi(s))}
\\
&\hspace{3.5cm}
\circ_{s}^{1\left\llbracket\mathbf{Pth}_{\boldsymbol{\mathcal{C}}^{(2)}}^{\mathbf{g}^{(2)}\circ\mathbf{f}^{(2)}}\right\rrbracket}
\left\llbracket
\left(
g^{(2)\flat}_{\varphi} \circ f^{(2)}
\right)_{s}^{\flat} \left(
\mathfrak{P}^{(2)0,i}
\right)
\right\rrbracket_{\psi(\varphi(s))}
\tag{2}
\\
&=
\left\llbracket
\left(g^{(2)\flat}_{\varphi} \circ f^{(2)\flat}\right)_{s} \left(
\mathfrak{P}^{(2)i+1, \bb{\mathfrak{P}^{(2)}}-1}
\right)
\right\rrbracket_{\psi(\varphi(s))}
\\
&\hspace{3.5cm}
\circ_{s}^{1\left\llbracket\mathbf{Pth}_{\boldsymbol{\mathcal{C}}^{(1)}}^{\mathbf{g}^{(2)}\circ\mathbf{f}^{(2)}}\right\rrbracket}
\left[
\left(g^{(2)\flat}_{\varphi} \circ f^{(2)\flat}\right)_{s} \left(
\mathfrak{P}^{(2)0,i}
\right)
\right\rrbracket_{\psi(\varphi(s))}
\tag{3}
\\
&=
\left\llbracket
g_{\varphi(s)}^{(2)\flat} \left(
f_{s}^{(2)\flat} \left(
\mathfrak{P}^{(2)i+1, \bb{\mathfrak{P}^{(2)}}-1}
\right)
\right)
\right\rrbracket_{\psi(\varphi(s))}
\\
&\hspace{3.5cm}
\circ_{\varphi(s)}^{1\llbracket\mathbf{Pth}_{\boldsymbol{\mathcal{C}}^{(2)}}^{\mathbf{g}^{(2)}}\rrbracket}
\left\llbracket
g_{\varphi(s)}^{(2)\flat} \left(
f_{s}^{(2)\flat} \left(
\mathfrak{P}^{(2)0,i}
\right)
\right)
\right\rrbracket_{\psi(\varphi(s))}
\tag{4}
\\
&=
\left\llbracket
g_{\varphi(s)}^{(2)\flat} \left(
f_{s}^{(2)\flat} \left(
\mathfrak{P}^{(2)i+1, \bb{\mathfrak{P}^{(2)}}-1}
\right)
\circ_{\varphi(s)}^{1\mathbf{Pth}_{\boldsymbol{\mathcal{B}}^{(2)}}}
f_{s}^{(2)\flat} \left(
\mathfrak{P}^{(2)0,i}
\right)
\right)
\right\rrbracket_{\psi(\varphi(s))}
\tag{5}
\\
&=
\left\llbracket
g_{\varphi(s)}^{(2)\flat} \left(
f_{s}^{(2)\flat} \left(
\mathfrak{P}^{(2)i+1, \bb{\mathfrak{P}^{(2)}}-1}
\circ_{s}^{1\mathbf{Pth}_{\boldsymbol{\mathcal{A}}^{(2)}}}
\mathfrak{P}^{(2)0,i}
\right)
\right)
\right\rrbracket_{\psi(\varphi(s))}
\tag{6}
\\
&=
\left\llbracket
g_{\varphi(s)}^{(2)\flat} \left(
f_{s}^{(2)\flat} \left(
\mathfrak{P}^{(2)}
\right)
\right)
\right\rrbracket_{\psi(\varphi(s))}.
\tag{7}
\end{align*}
\end{flushleft}

The first equality follows from the fact that $\mathfrak{P}^{(2)} = \mathfrak{P}^{(2)i+1, \bb{\mathfrak{P}^{(2)}}-1} \circ_{s}^{1\mathbf{Pth}_{\boldsymbol{\mathcal{A}}^{(2)}}} \mathfrak{P}^{(2)0,i}$;
the second equality follows from Proposition~\ref{PHomPthExtKer};
the third equality follows by Artinian induction;
the fourth equality unravels the definition of the $s$-th component of the composition of $S$-sorted mappings and note that, for every sort $s$ in $S$,
$$
\circ_{s}^{1\llbracket\mathbf{Pth}_{\boldsymbol{\mathcal{C}}^{(2)}}^{\mathbf{g}^{(2)}\circ\mathbf{f}^{(2)}}\rrbracket}
=
\circ_{\psi(\varphi(s))}^{1\llbracket\mathbf{Pth}_{\boldsymbol{\mathcal{C}}^{(2)}}\rrbracket}
=
\circ_{\varphi(s)}^{1\llbracket\mathbf{Pth}_{\boldsymbol{\mathcal{C}}^{(2)}}^{\mathbf{g}^{(2)}}\rrbracket};
$$
the fifth equality follows from Proposition~\ref{PHomPthExtKer};
the sixth equality recovers the definition of $f^{(2)\flat}$;
finally, the last equality recovers the definition of $\mathfrak{P}^{(2)}$;

Case (2.1) follows.

Therefore we are left with the case of $\mathfrak{P}^{(2)}$ being a head-constant echelonless second-order path. It could be the case that (2.2) $\mathfrak{P}^{(2)}$ is not coherent. Then let $i\in\bb{\mathfrak{P}^{(2)}}$ be the maximum index for which the subpath $\mathfrak{P}^{(2)0,i}$ of $\mathfrak{P}^{(2)}$ is a coherent head-constant echelonless second-order path. Note that the pairs $(\mathfrak{P}^{(2),0,i}, s)$ and $(\mathfrak{P}^{(2), i+1, \bb{\mathfrak{P}^{(2)}}-1}, s)$ $\prec_{\mathbf{Pth}_{\boldsymbol{\mathcal{A}}^{(2)}}}$-precede the pair $(\mathfrak{P}^{(2)},s)$. Therefore, by induction, the equalities
$$
\left\llbracket
\left(g^{(2)\flat}_{\varphi} \circ f^{(2)}\right)^{\flat}_{s}\left(
\mathfrak{P}^{(2)0,i}
\right)
\right\rrbracket_{\psi(\varphi(s))}
=
\left\llbracket
g^{(2)\flat}_{\varphi(s)} \circ f^{(2)\flat}_{s}\left(
\mathfrak{P}^{(2)0,i}
\right)
\right\rrbracket_{\psi(\varphi(s))}
$$
$$
\left\llbracket
\left(g^{(2)\flat}_{\varphi} \circ f^{(2)}\right)^{\flat}_{s}\left(
\mathfrak{P}^{(2)i+1,\bb{\mathfrak{P}^{(2)}}-1}
\right)
\right\rrbracket_{\psi(\varphi(s))}
=
\left\llbracket
g^{(2)\flat}_{\varphi(s)} \circ f^{(2)\flat}_{s}\left(
\mathfrak{P}^{(2)i+1,\bb{\mathfrak{P}^{(2)}}-1}
\right)
\right\rrbracket_{\psi(\varphi(s))}
$$
hold. Then, the following chain of equalities holds
\begin{flushleft}
$
\left\llbracket
\left(
g^{(2)\flat}_{\varphi} \circ f^{(2)}
\right)_{s}^{\flat} \left(
\mathfrak{P}^{(2)}
\right)
\right\rrbracket_{\psi(\varphi(s))}
$
\allowdisplaybreaks
\begin{align*}
&=
\left\llbracket
\left(
g^{(2)\flat}_{\varphi} \circ f^{(2)}
\right)_{s}^{\flat} \left(
\mathfrak{P}^{(2)i+1, \bb{\mathfrak{P}^{(2)}}-1} \circ_{s}^{1\mathbf{Pth}_{\boldsymbol{\mathcal{A}}^{(2)}}} \mathfrak{P}^{(2)0,i}
\right)
\right\rrbracket_{\psi(\varphi(s))}
\tag{1}
\\
&=
\left\llbracket
\left(
g^{(2)\flat}_{\varphi} \circ f^{(2)}
\right)_{s}^{\flat} \left(
\mathfrak{P}^{(2)i+1, \bb{\mathfrak{P}^{(2)}}-1}
\right)
\right\rrbracket_{\psi(\varphi(s))}
\\
&\hspace{3.5cm}
\circ_{s}^{1\left\llbracket\mathbf{Pth}_{\boldsymbol{\mathcal{C}}^{(2)}}^{\mathbf{g}^{(2)}\circ\mathbf{f}^{(2)}}\right\rrbracket}
\left\llbracket
\left(
g^{(2)\flat}_{\varphi} \circ f^{(2)}
\right)_{s}^{\flat} \left(
\mathfrak{P}^{(2)0,i}
\right)
\right\rrbracket_{\psi(\varphi(s))}
\tag{2}
\\
&=
\left\llbracket
\left(g^{(2)\flat}_{\varphi} \circ f^{(2)\flat}\right)_{s} \left(
\mathfrak{P}^{(2)i+1, \bb{\mathfrak{P}^{(2)}}-1}
\right)
\right\rrbracket_{\psi(\varphi(s))}
\\
&\hspace{3.5cm}
\circ_{s}^{1\left\llbracket\mathbf{Pth}_{\boldsymbol{\mathcal{C}}^{(1)}}^{\mathbf{g}^{(2)}\circ\mathbf{f}^{(2)}}\right\rrbracket}
\left[
\left(g^{(2)\flat}_{\varphi} \circ f^{(2)\flat}\right)_{s} \left(
\mathfrak{P}^{(2)0,i}
\right)
\right\rrbracket_{\psi(\varphi(s))}
\tag{3}
\\
&=
\left\llbracket
g_{\varphi(s)}^{(2)\flat} \left(
f_{s}^{(2)\flat} \left(
\mathfrak{P}^{(2)i+1, \bb{\mathfrak{P}^{(2)}}-1}
\right)
\right)
\right\rrbracket_{\psi(\varphi(s))}
\\
&\hspace{3.5cm}
\circ_{\varphi(s)}^{1\llbracket\mathbf{Pth}_{\boldsymbol{\mathcal{C}}^{(2)}}^{\mathbf{g}^{(2)}}\rrbracket}
\left\llbracket
g_{\varphi(s)}^{(2)\flat} \left(
f_{s}^{(2)\flat} \left(
\mathfrak{P}^{(2)0,i}
\right)
\right)
\right\rrbracket_{\psi(\varphi(s))}
\tag{4}
\\
&=
\left\llbracket
g_{\varphi(s)}^{(2)\flat} \left(
f_{s}^{(2)\flat} \left(
\mathfrak{P}^{(2)i+1, \bb{\mathfrak{P}^{(2)}}-1}
\right)
\circ_{\varphi(s)}^{1\mathbf{Pth}_{\boldsymbol{\mathcal{B}}^{(2)}}}
f_{s}^{(2)\flat} \left(
\mathfrak{P}^{(2)0,i}
\right)
\right)
\right\rrbracket_{\psi(\varphi(s))}
\tag{5}
\\
&=
\left\llbracket
g_{\varphi(s)}^{(2)\flat} \left(
f_{s}^{(2)\flat} \left(
\mathfrak{P}^{(2)i+1, \bb{\mathfrak{P}^{(2)}}-1}
\circ_{s}^{1\mathbf{Pth}_{\boldsymbol{\mathcal{A}}^{(2)}}}
\mathfrak{P}^{(2)0,i}
\right)
\right)
\right\rrbracket_{\psi(\varphi(s))}
\tag{6}
\\
&=
\left\llbracket
g_{\varphi(s)}^{(2)\flat} \left(
f_{s}^{(2)\flat} \left(
\mathfrak{P}^{(2)}
\right)
\right)
\right\rrbracket_{\psi(\varphi(s))}.
\tag{7}
\end{align*}
\end{flushleft}

The first equality follows from the fact that $\mathfrak{P}^{(2)} = \mathfrak{P}^{(2)i+1, \bb{\mathfrak{P}^{(2)}}-1} \circ_{s}^{1\mathbf{Pth}_{\boldsymbol{\mathcal{A}}^{(2)}}} \mathfrak{P}^{(2)0,i}$;
the second equality follows from Proposition~\ref{PDHomPthExtKer};
the third equality follows by Artinian induction;
the fourth equality unravels the definition of the $s$-th component of the composition of $S$-sorted mappings and noting that, for every sort $s$ in $S$,
$$
\circ_{s}^{0\llbracket\mathbf{Pth}_{\boldsymbol{\mathcal{C}}^{(2)}}^{\mathbf{g}^{(2)}\circ\mathbf{f}^{(2)}}\rrbracket}
=
\circ_{\psi(\varphi(s))}^{\llbracket\mathbf{Pth}_{\boldsymbol{\mathcal{C}}^{(2)}}\rrbracket}
=
\circ_{\varphi(s)}^{0\llbracket\mathbf{Pth}_{\boldsymbol{\mathcal{C}}^{(2)}}^{\mathbf{g}^{(2)}}\rrbracket};
$$
the fifth equality follows from Proposition~\ref{PHomPthExtKer};
the sixth equality recovers the definition of $f^{(2)\flat}$;
finally, the last equality recovers the definition of $\mathfrak{P}^{(2)}$;

Case (2.2) follows.

Therefore we are left with the case (2.3) of $\mathfrak{P}^{(2)}$ being a coherent head-constant echelonless second-order path. Under this setting, the conditions for the second-order extraction algorithm, that is, Lemma~\ref{LDPthExtract}, are fulfilled. Then there exists a unique word $\mathbf{s} \in S^{\star} - \{\lambda\}$ and a unique operation symbol $\tau \in \Sigma_{\mathbf{s}, s}^{\boldsymbol{\mathcal{A}}^{(1)}}$ associated to $\mathfrak{P}^{(2)}$. Let $(\mathfrak{P}^{(2)}_{j})_{j \in \bb{\mathbf{s}}}$ be the family of second-order paths in $\mathrm{Pth}_{\boldsymbol{\mathcal{A}}^{(2)}, \mathbf{s}}$, which, in virtue of Lemma~\ref{LDPthExtract}, we can extract from $\mathfrak{P}^{(2)}$. Note that, for every $j \in \bb{\mathbf{s}}$, the pair $(\mathfrak{P}_{j}^{(2)}, s_{j})$ $\prec_{\mathbf{Pth}_{\boldsymbol{\mathcal{A}}^{(2)}}}$- precede the pair $(\mathfrak{P}^{(2)}, s)$. Therefore, by induction, for every $j \in \bb{\mathbf{s}}$, the equality
$$
\left\llbracket
\left(g^{(2)\flat}_{\varphi} \circ f^{(2)}\right)^{\flat}_{s_{j}}\left(
\mathfrak{P}^{(2)}_{j}
\right)
\right\rrbracket_{\psi(\varphi(s_{j}))}
=
\left\llbracket
g^{(2)\flat}_{\varphi(s_{j})} \circ f^{(2)\flat}_{s_{j}}\left(
\mathfrak{P}^{(2)}_{j}
\right)
\right\rrbracket_{\psi(\varphi(s_{j}))}
$$
hold. Then, we consider different cases for $\tau$ according to the cases presented in Definition~\ref{DUTrans}.

If $\tau = \sigma$ is an operation symbol in $\Sigma_{\mathbf{s}, s}$, then the following chain of equalities holds
\begin{flushleft}
$
\left\llbracket
\left(
g^{(2)\flat}_{\varphi} \circ f^{(2)}
\right)_{s}^{\flat} \left(
\mathfrak{P}^{(2)}
\right)
\right\rrbracket_{\psi(\varphi(s))}
$
\allowdisplaybreaks
\begin{align*}
&=
\left\llbracket
\sigma^{\mathbf{Pth}_{\boldsymbol{\mathcal{C}}^{(2)}}^{\mathbf{g}^{(2)} \circ \mathbf{f}^{(2)}}} \left(
\left(
\left(
g_{\varphi}^{(2)\flat} \circ
f^{(2)}
\right)_{s_{j}}^{\flat}
\left(
\mathfrak{P}^{(2)}_{j}
\right)
\right)_{j \in \bb{\mathbf{s}}}
\right)
\right\rrbracket_{\psi(\varphi(s))}
\tag{1}
\\
&=
\sigma^{\llbracket\mathbf{Pth}_{\boldsymbol{\mathcal{C}}^{(2)}}^{\mathbf{g}^{(2)} \circ \mathbf{f}^{(2)}}\rrbracket} \left(
\left(
\left\llbracket
\left(
g_{\varphi}^{(2)\flat} \circ
f^{(2)}
\right)_{s_{j}}^{\flat}
\left(
\mathfrak{P}^{(2)}_{j}
\right)
\right\rrbracket_{\psi(\varphi(s_{j}))}
\right)_{j \in \bb{\mathbf{s}}}
\right)
\tag{2}
\\
&=
\sigma^{\llbracket\mathbf{Pth}_{\boldsymbol{\mathcal{C}}^{(2)}}^{\mathbf{g}^{(2)} \circ \mathbf{f}^{(2)}}\rrbracket} \left(
\left(
\left\llbracket
g^{(2)\flat}_{\varphi(s_{j})} \circ 
f^{(2)\flat}_{s_{j}}\left(
\mathfrak{P}^{(2)}_{j}
\right)
\right\rrbracket_{\psi(\varphi(s_{j}))}
\right)_{j \in \bb{\mathbf{s}}}
\right)
\tag{3}
\\
&=
\left\llbracket
\sigma^{\mathbf{Pth}_{\boldsymbol{\mathcal{C}}^{(2)}}^{\mathbf{g}^{(2)} \circ \mathbf{f}^{(2)}}} \left(
\left(
g_{\varphi(s_{j)}}^{(2)\flat} \circ
f^{(2)\flat}_{s_{j}}
\left(
\mathfrak{P}^{(2)}_{j}
\right)
\right)_{j \in \bb{\mathbf{s}}}
\right)
\right\rrbracket_{\psi(\varphi(s))}
\tag{4}
\\
&=
\left\llbracket
\sigma^{\mathbf{d}_{\mathfrak{d}}^{\ast}(\mathbf{c}_{\mathfrak{d}}^{\ast}(\mathbf{Pth}_{\boldsymbol{\mathcal{C}}^{(2)}}))} \left(
\left(
g_{\varphi(s_{j})}^{(2)\flat} \circ
f^{(2)\flat}_{s_{j}}
\left(
\mathfrak{P}^{(2)}_{j}
\right)
\right)_{j \in \bb{\mathbf{s}}}
\right)
\right\rrbracket_{\psi(\varphi(s))}
\tag{5}
\\
&=
\left\llbracket
g_{\varphi(s)}^{(2)\flat} \left(
\sigma^{\mathbf{c}_{\mathfrak{d}}^{\ast}(\mathbf{Pth}_{\boldsymbol{\mathcal{B}}^{(2)}})} \left(
\left(
f^{(2)\flat}_{s_{j}}
\left(
\mathfrak{P}^{(2)}_{j}
\right)
\right)_{j \in \bb{\mathbf{s}}}
\right)
\right)
\right\rrbracket_{\psi(\varphi(s))}
\tag{6}
\\
&=
\left\llbracket
g^{(2)\flat}_{\varphi(s)} \left( 
f^{(2)\flat}_{s}\left(
\mathfrak{P}^{(2)}
\right)
\right)
\right\rrbracket_{\psi(\varphi(s))}
\tag{7}
\end{align*}
\end{flushleft}

The first equality unravels the definition of the second-order path extension mapping introduced in Proposition~\ref{PDPthExt} considering the second-order morphism $\mathbf{g}^{(2)} \circ \mathbf{f}^{(2)}$;
the second equality recovers the definition of the interpretation of $\sigma$ in the partial $\Sigma^{\boldsymbol{\mathcal{A}}^{(2)}}$-algebra $\llbracket \mathbf{Pth}_{\boldsymbol{\mathcal{C}}^{(2)}}^{\mathbf{g}^{(2)} \circ \mathbf{f}^{(2)}} \rrbracket$ introduced in Proposition~\ref{PDQPthBDCatAlg};
the third equality follows by Artinian induction;
the fourth equality unravels the definition of the interpretation of $\sigma$ in the partial $\Sigma^{\boldsymbol{\mathcal{A}}^{(2)}}$-algebra $\llbracket \mathbf{Pth}_{\boldsymbol{\mathcal{C}}^{(2)}}^{\mathbf{g}^{(2)} \circ \mathbf{f}^{(2)}} \rrbracket$ introduced in Proposition~\ref{PDQPthBDCatAlg};
the fifth equality unravels the definition of the interpretation of $\sigma$ in the partial $\Sigma^{\boldsymbol{\mathcal{A}}^{(2)}}$-algebra $\mathbf{Pth}_{\boldsymbol{\mathcal{C}}^{(2)}}^{\mathbf{g}^{(2)} \circ \mathbf{f}^{(2)}}$ introduced in Proposition~\ref{PDPthBDCatAlg};
the sixth equality follows from the fact that, according to Proposition~\ref{PDPthExtHom}, $g^{(2)\flat}$ is a $\Lambda$-homomorphism, thus, according to Proposition~\ref{PFunSig}, $\mathbf{c}^{\ast}_{\mathfrak{d}}(g^{(2)\flat})$ is a $\Sigma$-homomorphism;
finally, the last equality follows from the fact that, according to Proposition~\ref{PDPthExtHom}, $f^{(2)\flat}$ is a $\Sigma$-homomorphism.

If $\tau = \circ^{0}_{s}$ is the operation symbol of $0$-composition, then the following chain of equalities holds
\begin{flushleft}
$
\left\llbracket
\left(
g^{(2)\flat}_{\varphi} \circ f^{(2)}
\right)_{s}^{\flat} \left(
\mathfrak{P}^{(2)}
\right)
\right\rrbracket_{\psi(\varphi(s))}
$
\allowdisplaybreaks
\begin{align*}
&=
\left\llbracket
\left(
g^{(2)\flat}_{\varphi} \circ f^{(2)}
\right)_{s}^{\flat} \left(
\mathfrak{P}^{(2)}_{1} \circ_{s}^{0\mathbf{Pth}_{\boldsymbol{\mathcal{A}}^{(2)}}} \mathfrak{P}^{(2)}_{0}
\right)
\right\rrbracket_{\psi(\varphi(s))}
\tag{1}
\\
&=
\left\llbracket
\left(
g^{(2)\flat}_{\varphi} \circ f^{(2)}
\right)_{s}^{\flat} \left(
\mathfrak{P}^{(2)}_{1}
\right)
\right\rrbracket_{\psi(\varphi(s))}
\\
&\hspace{3.5cm}
\circ_{s}^{0\left\llbracket\mathbf{Pth}_{\boldsymbol{\mathcal{C}}^{(2)}}^{\mathbf{g}^{(2)}\circ\mathbf{f}^{(2)}}\right\rrbracket}
\left\llbracket
\left(
g^{(2)\flat}_{\varphi} \circ f^{(2)}
\right)_{s}^{\flat} \left(
\mathfrak{P}^{(2)}_{0}
\right)
\right\rrbracket_{\psi(\varphi(s))}
\tag{2}
\\
&=
\left\llbracket
\left(g^{(2)\flat}_{\varphi} \circ f^{(2)\flat}\right)_{s} \left(
\mathfrak{P}^{(2)}_{1}
\right)
\right\rrbracket_{\psi(\varphi(s))}
\\
&\hspace{3.5cm}
\circ_{s}^{0\left\llbracket\mathbf{Pth}_{\boldsymbol{\mathcal{C}}^{(1)}}^{\mathbf{g}^{(2)}\circ\mathbf{f}^{(2)}}\right\rrbracket}
\left[
\left(g^{(2)\flat}_{\varphi} \circ f^{(2)\flat}\right)_{s} \left(
\mathfrak{P}^{(2)}_{0}
\right)
\right\rrbracket_{\psi(\varphi(s))}
\tag{3}
\\
&=
\left\llbracket
g_{\varphi(s)}^{(2)\flat} \left(
f_{s}^{(2)\flat} \left(
\mathfrak{P}^{(2)}_{1}
\right)
\right)
\right\rrbracket_{\psi(\varphi(s))}
\\
&\hspace{3.5cm}
\circ_{\varphi(s)}^{0\llbracket\mathbf{Pth}_{\boldsymbol{\mathcal{C}}^{(2)}}^{\mathbf{g}^{(2)}}\rrbracket}
\left\llbracket
g_{\varphi(s)}^{(2)\flat} \left(
f_{s}^{(2)\flat} \left(
\mathfrak{P}^{(2)}_{0}
\right)
\right)
\right\rrbracket_{\psi(\varphi(s))}
\tag{4}
\\
&=
\left\llbracket
g_{\varphi(s)}^{(2)\flat} \left(
f_{s}^{(2)\flat} \left(
\mathfrak{P}^{(2)}_{1}
\right)
\circ_{\varphi(s)}^{0\mathbf{Pth}_{\boldsymbol{\mathcal{B}}^{(2)}}}
f_{s}^{(2)\flat} \left(
\mathfrak{P}^{(2)}_{0}
\right)
\right)
\right\rrbracket_{\psi(\varphi(s))}
\tag{5}
\\
&=
\left\llbracket
g_{\varphi(s)}^{(2)\flat} \left(
f_{s}^{(2)\flat} \left(
\mathfrak{P}^{(2)}_{1}
\circ_{s}^{0\mathbf{Pth}_{\boldsymbol{\mathcal{A}}^{(2)}}}
\mathfrak{P}^{(2)}_{0}
\right)
\right)
\right\rrbracket_{\psi(\varphi(s))}
\tag{6}
\\
&=
\left\llbracket
g_{\varphi(s)}^{(2)\flat} \left(
f_{s}^{(2)\flat} \left(
\mathfrak{P}^{(2)}
\right)
\right)
\right\rrbracket_{\psi(\varphi(s))}.
\tag{7}
\end{align*}
\end{flushleft}

The first equality unravels the definition of $\mathfrak{P}^{(2)}$;
the second equality follows from Proposition~\ref{PDHomPthExtKer};
the third equality follows by Artinian induction;
the fourth equality unravels the definition of the $s$-th component of the composition of $S$-sorted mappings and noting that, for every sort $s$ in $S$,
$$
\circ_{s}^{1\llbracket\mathbf{Pth}_{\boldsymbol{\mathcal{C}}^{(2)}}^{\mathbf{g}^{(2)}\circ\mathbf{f}^{(2)}}\rrbracket}
=
\circ_{\psi(\varphi(s))}^{\llbracket\mathbf{Pth}_{\boldsymbol{\mathcal{C}}^{(2)}}\rrbracket}
=
\circ_{\varphi(s)}^{1\llbracket\mathbf{Pth}_{\boldsymbol{\mathcal{C}}^{(2)}}^{\mathbf{g}^{(2)}}\rrbracket};
$$
the fifth equality follows from Proposition~\ref{PDHomPthExtKer};
the sixth equality recovers the definition of $f^{(2)\flat}$;
finally, the last equality recovers the definition of $\mathfrak{P}^{(2)}$;

Case (2.3) follows.

This completes case (2).

This completes the proof.
\end{proof}

In the next corollary, we give equivalent expressions for the homomorphic path extension mapping of a second-order composition of second-order morphisms in terms of the respective homomorphic path extension mappings of the composites.

\begin{corollary}\label{CDQPthExtComp}
Let $\mathbf{f}^{(2)}=(\varphi, c, (f^{(i)})_{i\in 3})$ be a second-order morphism from $\boldsymbol{\mathcal{A}}^{(2)}$ to $\boldsymbol{\mathcal{B}}^{(2)}$ and $\mathbf{g}^{(2)}=(\psi, d, (g^{(i)})_{i\in 3})$ a second-order morphism from $\boldsymbol{\mathcal{B}}^{(2)}$ to $\boldsymbol{\mathcal{C}}^{(2)}$. Then the following three mappings are equal:
\begin{enumerate}
\item
$\mathrm{pr}^{\llbracket\cdot\rrbracket}_{\boldsymbol{\mathcal{C}}^{(2)}, \psi\circ\varphi} \circ \left(g^{(2)\flat}_{\varphi} \circ f^{(2)}\right)^{\flat}$;
\item
$\left(\mathrm{pr}^{\llbracket\cdot\rrbracket}_{\boldsymbol{\mathcal{C}}^{(2)}, \psi} \circ g^{(2)\flat}\right)_{\varphi} \circ f^{(2)\flat}$;
\item
$g^{(2)@}_{\varphi} \circ \left(\mathrm{pr}^{\llbracket\cdot\rrbracket}_{\boldsymbol{\mathcal{B}}^{(2)}, \varphi} \circ f^{(2)\flat}\right)$.
\end{enumerate}
\end{corollary}

Finally, in the next corollary we give an expression for the quotient path extension mapping of the second-order composition of second-order morphisms in terms of the respective quotient path extension mappings of the composites.

\begin{corollary}\label{CDQPthExtComp}
Let $\mathbf{f}^{(2)}=(\varphi, c, (f^{(i)})_{i\in 3})$ be a second-order morphism from $\boldsymbol{\mathcal{A}}^{(2)}$ to $\boldsymbol{\mathcal{B}}^{(2)}$ and $\mathbf{g}^{(2)}=(\psi, d, (g^{(i)})_{i\in 3})$ a second-order morphism from $\boldsymbol{\mathcal{B}}^{(2)}$ to $\boldsymbol{\mathcal{C}}^{(2)}$. Then 
$$
\left(
g^{(2)\flat}_{\varphi} \circ f^{(2)}
\right)^{@}
=
g^{\llbracket2\rrbracket@}_{\varphi} \circ f^{\llbracket2\rrbracket@}.
$$
\end{corollary}

We are now in a position where we can prove that the second-order rewriting systems and equivalence classes of second-order morphisms constitute a category. Before doing so, we show that the composition of classes of second-order morphisms is a well-defined binary operation and it does not depend on the choice of the representatives.

\begin{proposition}
\label{PCompRws2}
Let $\mathbf{f}^{(2)}$ and $\mathbf{f}'^{(2)}$ be second-order morphisms from $\boldsymbol{\mathcal{A}}^{(2)}$ to $\boldsymbol{\mathcal{B}}^{(2)}$ and $\mathbf{g}^{(2)}$ and $\mathbf{g}'^{(2)}$ second-order morphisms from $\boldsymbol{\mathcal{B}}^{(2)}$ to $\boldsymbol{\mathcal{C}}^{(2)}$. If $\mathbf{f}^{(2)} \cong^{(2)} \mathbf{f}'^{(2)}$ and $\mathbf{g}^{(2)} \cong^{(2)} \mathbf{g}'^{(2)}$, then
$
\mathbf{g}^{(2)} \circ \mathbf{f}^{(2)}
\cong^{(2)}
\mathbf{g}'^{(2)} \circ \mathbf{f}'^{(2)}.
$
Thus, the composition of their morphism classes $\llbracket\mathbf{g}^{(2)}\rrbracket \circ \llbracket\mathbf{f}^{(2)}\rrbracket = \llbracket\mathbf{g}^{(2)} \circ \mathbf{f}^{(2)}\rrbracket$ is well-defined and does not depend on the representatives of the second-order morphism classes.
\end{proposition}

\begin{proof}
Let $\mathbf{f}^{(2)} = (\mathbf{f}^{(1)}, f^{(2)})$ and $\mathbf{f}'^{(2)} = (\mathbf{f}'^{(1)}, f'^{(2)})$. Since $\mathbf{f}^{(2)}$ and $\mathbf{f}'^{(2)}$ are second-order equivalent, following Definition~\ref{DDMorEqv}, $\mathbf{f}^{(1)} \cong^{(1)} \mathbf{f}'^{(1)}$ and 
\begin{equation}
\mathrm{pr}^{\llbracket\cdot\rrbracket}_{\boldsymbol{\mathcal{B}}^{(2)}, \varphi} \circ
f^{(2)\flat}
=
\mathrm{pr}^{\llbracket\cdot\rrbracket}_{\boldsymbol{\mathcal{B}}^{(2)}, \varphi} \circ
f'^{(2)\flat}.
\tag{F}
\end{equation}
Moreover, let $\mathbf{g}^{(2)} = (\mathbf{g}^{(1)}, g^{(2)})$ and $\mathbf{g}'^{(2)} = (\mathbf{g}'^{(1)}, g'^{(2)})$. Since $\mathbf{g}^{(2)}$ and $\mathbf{g}'^{(2)}$ are second-order equivalent, following Definition~\ref{DDMorEqv}, $\mathbf{g}^{(1)} \cong^{(1)} \mathbf{g}'^{(1)}$ and 
\begin{equation}
\mathrm{pr}^{\llbracket\cdot\rrbracket}_{\boldsymbol{\mathcal{C}}^{(2)}, \psi}
\circ
g^{(2)\flat}
=
\mathrm{pr}^{\llbracket\cdot\rrbracket}_{\boldsymbol{\mathcal{C}}^{(2)}, \psi}
\circ
g'^{(2)\flat}.
\tag{G}
\end{equation}

Following Definition~\ref{DCompRws2}, the respective composite morphisms are given by
$$
\mathbf{g}^{(2)} \circ \mathbf{f}^{(2)}
=
\left(
\mathbf{g}^{(1)} \circ \mathbf{f}^{(1)},
g^{(2)\flat}_{\varphi} \circ f^{(2)}
\right)
$$
and
$$
\mathbf{g}'^{(2)} \circ \mathbf{f}'^{(2)}
=
\left(
\mathbf{g}'^{(1)} \circ \mathbf{f}'^{(1)},
g'^{(2)\flat}_{\varphi} \circ f'^{(2)}
\right).
$$

Recall that, according to Proposition~\ref{PCompRws1}, it follows that $\mathbf{g}^{(1)} \circ \mathbf{f}^{(1)}$ and  $\mathbf{g}'^{(1)} \circ \mathbf{f}'^{(1)}$ are first-order equivalent. Therefore, by Definition~\ref{DDMorEqv}, to show that $\mathbf{g}^{(2)} \circ \mathbf{f}^{(2)}$ and $\mathbf{g}'^{(2)} \circ \mathbf{f}'^{(2)}$ are second-order equivalent, it suffices to show that
$$
\mathrm{pr}^{\llbracket\cdot\rrbracket}_{\boldsymbol{\mathcal{C}}^{(2)}, \psi \circ \varphi}
\circ
\left(
g^{(2)\flat}_{\varphi} \circ 
f^{(2)}
\right)^{\flat}
=
\mathrm{pr}^{\llbracket\cdot\rrbracket}_{\boldsymbol{\mathcal{C}}^{(2)}, \psi \circ \varphi}
\circ
\left(
g'^{(2)\flat}_{\varphi} \circ 
f'^{(2)}
\right)^{\flat}
$$

The following chain of equalities holds
\allowdisplaybreaks
\begin{align*}
\mathrm{pr}^{\llbracket\cdot\rrbracket}_{\boldsymbol{\mathcal{C}}^{(2)}, \psi \circ \varphi}
\circ
\left(
g^{(2)\flat}_{\varphi} \circ 
f^{(2)}
\right)^{\flat}
&=
\mathrm{pr}^{\llbracket\cdot\rrbracket}_{\boldsymbol{\mathcal{C}}^{(2)}, \psi \circ \varphi}
\circ
g^{(2)\flat}_{\varphi} \circ
f^{(2)\flat}
\tag{1}
\\
&=
\left(
\mathrm{pr}^{\llbracket\cdot\rrbracket}_{\boldsymbol{\mathcal{C}}^{(2)}, \psi}
\right)_{\varphi} \circ
g^{(2)\flat}_{\varphi} \circ
f^{(2)\flat}
\tag{2}
\\
&=
\left(
\mathrm{pr}^{\llbracket\cdot\rrbracket}_{\boldsymbol{\mathcal{C}}^{(2)}, \psi}
\circ
g^{(2)\flat}\right)_{\varphi} \circ
f^{(2)\flat}
\tag{3}
\\
&=
\left(
\mathrm{pr}^{\llbracket\cdot\rrbracket}_{\boldsymbol{\mathcal{C}}^{(2)}, \psi}
\circ
g'^{(2)\flat}\right)_{\varphi} \circ
f^{(2)\flat}
\tag{4}
\\
&=
\left(
g'^{\llbracket 2 \rrbracket@}
\circ
\mathrm{pr}^{\llbracket\cdot\rrbracket}_{\boldsymbol{\mathcal{B}}^{(2)}}\right)_{\varphi} \circ
f^{(2)\flat}
\tag{5}
\\
&=
\left(
g'^{\llbracket 2 \rrbracket@}
\right)_{\varphi}
\circ
\mathrm{pr}^{\llbracket\cdot\rrbracket}_{\boldsymbol{\mathcal{B}}^{(2)}, \varphi} \circ
f^{(2)\flat}
\tag{6}
\\
&=
\left(
g'^{\llbracket 2 \rrbracket@}
\right)_{\varphi}
\circ
\mathrm{pr}^{\llbracket\cdot\rrbracket}_{\boldsymbol{\mathcal{B}}^{(2)}, \varphi} \circ
f'^{(2)\flat}
\tag{7}
\\
&=
\left(
g'^{\llbracket 2 \rrbracket@}
\circ
\mathrm{pr}^{\llbracket\cdot\rrbracket}_{\boldsymbol{\mathcal{B}}^{(2)}} \right)_{\varphi} \circ
f'^{(2)\flat}
\tag{8}
\\
&=
\left(
\mathrm{pr}^{\llbracket\cdot\rrbracket}_{\boldsymbol{\mathcal{C}}^{(2)}, \psi}
\circ
g'^{(2)\flat}  \right)_{\varphi}
\circ
f'^{(2)\flat}
\tag{9}
\\
&=
\left(
\mathrm{pr}^{\llbracket\cdot\rrbracket}_{\boldsymbol{\mathcal{C}}^{(2)}, \psi}
\right)_{\varphi}
\circ
g'^{(2)\flat}_{\varphi}
\circ
f'^{(2)\flat}
\tag{10}
\\
&=
\mathrm{pr}^{\llbracket\cdot\rrbracket}_{\boldsymbol{\mathcal{C}}^{(2)}, \psi \circ \phi}
\circ
g'^{(2)\flat}_{\varphi}
\circ
f'^{(2)\flat}
\tag{11}
\\
&=
\mathrm{pr}^{\llbracket\cdot\rrbracket}_{\boldsymbol{\mathcal{C}}^{(2)}, \psi \circ \phi}
\circ
\left(
g'^{(2)\flat}_{\varphi}
\circ
f'^{(2)}
\right)^{\flat}.
\tag{12}
\end{align*}

The first equality follows from Proposition~\ref{PDPthExtComp};
the second equality follows from the fact that, according to Proposition~\ref{PMSetFunc}, $\mathrm{MSet}$ is a contravariant functor;
the third equality follows from the fact that, according to Proposition~\ref{PDeltaPhiFunc}, $\Delta_{\varphi}$ is a covariant functor;
the fourth equality follows from Equation~(G);
the fifth equality follows from the definition of $g'^{\llbracket 2 \rrbracket@}$, introduced in Definition~\ref{DDQPthExt};
the sixth equality follows from the fact that, according to Proposition~\ref{PDeltaPhiFunc}, $\Delta_{\varphi}$ is a covariant functor;
the seventh equality follows from Equation~(F);
the eighth equality follows from the fact that, according to Proposition~\ref{PDeltaPhiFunc}, $\Delta_{\varphi}$ is a covariant functor;
the ninth equality follows from the definition of $g'^{(2)@}$, introduced in Definition~\ref{DDQPthExt};  
the tenth equality follows from the fact that, according to Proposition~\ref{PDeltaPhiFunc}, $\Delta_{\varphi}$ is a covariant functor;
the eleventh equality follows from the fact that, according to Proposition~\ref{PMSetFunc}, $\mathrm{MSet}$ is a contravariant functor;
finally, the last equality follows from Proposition~\ref{PDPthExtComp}.

This completes the proof.
\end{proof}

We now show that second-order many-sorted rewriting systems and classes of second-order morphisms between them constitute a category.

\begin{proposition}\label{PRws2Cat}
The second-order many-sorted rewriting systems together with the classes of second-order morphisms between second-order many-sorted rewriting systems constitute a category, denoted by $\mathsf{Rws}^{\llbracket 2 \rrbracket}_{\mathfrak{d}}$.
\end{proposition}

\begin{proof}
That domains and codomains respect identities and compositions follows from the definitions of second-order identity morphism and second-order composition of second-order morphisms introduced in Definitions~\ref{DIdRws2} and \ref{DCompRws2}. Thus, all that remains to be proven is that the class of the second-order identity morphism at a second-order many-sorted rewriting system acts as a unit element and that the second-order composition of classes of second-order morphisms is associative.

\textsf{Unit element.}

Let $\mathbf{f}^{(2)}=(\varphi, c, (f^{(i)})_{i\in 3})$ be a second-order morphism from $\boldsymbol{\mathcal{A}}^{(2)}=(\boldsymbol{\mathcal{A}}^{(1)}, \mathcal{A}^{(2)})$ to $\boldsymbol{\mathcal{B}}^{(2)}=(\boldsymbol{\mathcal{B}}^{(1)}, \mathcal{B}^{(2)})$. We need to prove that
\begin{align*}
\left\llbracket
\mathbf{f}^{(2)}
\right\rrbracket 
\circ 
\left\llbracket
\mathrm{id}^{\boldsymbol{\mathcal{A}}^{(2)}}
\right\rrbracket
&=
\left\llbracket
\mathbf{f}^{(2)}
\right\rrbracket
&&\mbox{ and }&
\left\llbracket
\mathrm{id}^{\boldsymbol{\mathcal{B}}^{(2)}}
\right\rrbracket
\circ
\left\llbracket
\mathbf{f}^{(2)}
\right\rrbracket 
&=
\left\llbracket
\mathbf{f}^{(2)}
\right\rrbracket.
\tag{Id}
\end{align*}

We will only prove the left hand side of (Id). The right hand side of (Id) is done similarly. According to Definition~\ref{DIdRws2}, the second-order identity morphism at $\boldsymbol{\mathcal{A}}^{(2)}$ is given by $(\mathrm{id}^{\boldsymbol{\mathcal{A}}^{(1)}}, \mathrm{ech}^{(2,\mathcal{A}^{(2)})}_{\boldsymbol{\mathcal{A}}^{(2)}})$. Moreover, according to Definition~\ref{DCompRws2} and Proposition~\ref{PCompRws2}, the second-order composition morphism is
$$
\left\llbracket
\mathbf{f}^{(2)}
\right\rrbracket
\circ
\left\llbracket
\mathrm{id}^{\boldsymbol{\mathcal{A}}^{(2)}}
\right\rrbracket
=
\left\llbracket
\mathbf{f}^{(2)}
\circ
\mathrm{id}^{\boldsymbol{\mathcal{A}}^{(2)}}
\right\rrbracket
=
\left\llbracket
\left(
\mathbf{f}^{(1)} \circ \mathrm{id}^{\boldsymbol{\mathcal{A}}^{(1)}},
f^{(2)\flat}_{\mathrm{id}^{S}} \circ \mathrm{ech}^{(2,\mathcal{A}^{(2)})}_{\boldsymbol{\mathcal{A}}^{(2)}}
\right)
\right\rrbracket.
$$

According Proposition~\ref{PRws1Cat}, $\mathsf{Rws}_{\mathfrak{d}}^{[1]}$ is a category. Thus, $\mathbf{f}^{(1)} \circ \mathrm{id}^{\boldsymbol{\mathcal{A}}^{(1)}} \cong^{(1)} \mathbf{f}^{(1)}$. Moreover, the following chain of equalities holds
\allowdisplaybreaks
\begin{align*}
f^{(2)\flat}_{\mathrm{id}^{S}}
\circ
\mathrm{ech}^{(2,\mathcal{A}^{(2)})}_{\boldsymbol{\mathcal{A}}^{(2)}}
&=
f^{(2)\flat} \circ \mathrm{ech}^{(2,\mathcal{A}^{(2)})}_{\boldsymbol{\mathcal{A}}^{(2)}}
\tag{1}
\\
&=
f^{(2)}.
\tag{2}
\end{align*}

The first equality follows from the fact that, by Proposition~\ref{PMSetFunc}, $\mathrm{MSet}$ is a functor;
the second equality follows from Proposition~\ref{PDPthExt}.

Thus, Equation~(Id) follows.

This proves that the class of a second-order identity morphism at a second-order many-sorted rewriting system act as the unit element for second-order composition of classes of second-order morphisms.

\textsf{Associativity.}

Let $\boldsymbol{\mathcal{A}}^{(2)}$, $\boldsymbol{\mathcal{B}}^{(2)}$, $\boldsymbol{\mathcal{C}}^{(2)}$ and $\boldsymbol{\mathcal{D}}^{(2)}$ be four second-order many-sorted rewriting systems
and let $\mathbf{f}^{(2)}$, $\mathbf{g}^{(2)}$ and $\mathbf{h}^{(2)}$ be second-order morphisms of the form
$$
\mathbf{f}^{(2)} \colon \boldsymbol{\mathcal{A}}^{(2)} \mor \boldsymbol{\mathcal{B}}^{(2)}, \,
\mathbf{g}^{(2)} \colon \boldsymbol{\mathcal{B}}^{(2)} \mor \boldsymbol{\mathcal{C}}^{(2)}
\mbox{ and }
\mathbf{h}^{(2)} \colon \boldsymbol{\mathcal{C}}^{(2)} \mor \boldsymbol{\mathcal{D}}^{(2)}
$$
where $\mathbf{f}^{(2)}$, $\mathbf{g}^{(2)}$ and $\mathbf{h}^{(2)}$ stand for
\begin{align*}
\mathbf{f}^{(2)}=(\varphi, c, (f^{(i)})_{i\in 3}), \,
\mathbf{g}^{(2)}=(\psi, d, (g^{(i)})_{i\in 3})
\\
\mbox{ and }
\mathbf{h}^{(2)}=(\theta, e, (h^{(i)})_{i\in 3}),
\mbox{ respectively.}
\end{align*}
We need to prove that 
\begin{equation}
\left\llbracket\mathbf{h}^{(2)}\right\rrbracket
\circ
\left(
\left\llbracket\mathbf{g}^{(2)}\right\rrbracket
\circ
\left\llbracket\mathbf{f}^{(2)}\right\rrbracket
\right)
=
\left(
\left\llbracket\mathbf{h}^{(1)}\right\rrbracket
\circ
\left\llbracket\mathbf{g}^{(1)}\right\rrbracket
\right)
\circ
\left\llbracket\mathbf{f}^{(1)}\right\rrbracket.
\tag{Assoc}
\end{equation}

Following Definition~\ref{DCompRws2} and Proposition~\ref{PCompRws2}, the left hand side composition of Equation~(Assoc) is
\allowdisplaybreaks
\begin{align*}
\left\llbracket\mathbf{h}^{(2)}\right\rrbracket
\circ
\left(
\left\llbracket\mathbf{g}^{(2)}\right\rrbracket
\circ
\left\llbracket\mathbf{f}^{(2)}\right\rrbracket
\right)
&=
\left\llbracket
\mathbf{h}^{(2)}
\circ
\left( 
\mathbf{g}^{(2)}
\circ
\mathbf{f}^{(2)}
\right)
\right\rrbracket
\\
&=
\left\llbracket
\left(
\mathbf{h}^{(1)} \circ \left( \mathbf{g}^{(1)} \circ \mathbf{f}^{(1)} \right),
h^{(2)\flat}_{\psi\circ\varphi} \circ g^{(2)\flat}_{\varphi} \circ f^{(2)}
\right)
\right\rrbracket.
\end{align*}
Similarly, the right hand side composition of Equation~(Assoc) is
\allowdisplaybreaks
\begin{align*}
\left(
\left\llbracket\mathbf{h}^{(2)}\right\rrbracket
\circ
\left\llbracket\mathbf{g}^{(2)}\right\rrbracket
\right)
\circ
\left\llbracket\mathbf{f}^{(2)}\right\rrbracket
&=
\left\llbracket
\left(
\mathbf{h}^{(2)}
\circ
\mathbf{g}^{(2)}
\right)
\circ
\mathbf{f}^{(2)}
\right\rrbracket
\\
&=
\left\llbracket
\left(
\left(\mathbf{h}^{(1)} \circ \mathbf{g}^{(1)} \right) \circ \mathbf{f}^{(1)},
\left(h^{(2)\flat}_{\psi} \circ g^{(2)}\right)^{\flat}_{\varphi} \circ f^{(2)}
\right)
\right\rrbracket.
\end{align*}

According Proposition~\ref{PRws1Cat}, $\mathsf{Rws}_{\mathfrak{d}}^{[1]}$ is a category. Thus,
$$
\mathbf{h}^{(1)} \circ \left( \mathbf{g}^{(1)} \circ \mathbf{f}^{(1)} \right)
\cong^{(1)}
\left(\mathbf{h}^{(1)} \circ \mathbf{g}^{(1)} \right) \circ \mathbf{f}^{(1)}.
$$
Thus, according to Definition~\ref{DDMorEqv}, all what remains to be proven is that,
$$
\mathrm{pr}^{\llbracket\cdot\rrbracket}_{\boldsymbol{\mathcal{D}}^{(2)}, \theta \circ \psi \circ \varphi} 
\circ 
\left(
h^{(2)\flat}_{\psi \circ \varphi}
\circ
\left(
g^{(2)\flat}_{\varphi}
\circ
f^{(2)}
\right)
\right)^{\flat}
=
\mathrm{pr}^{\llbracket\cdot\rrbracket}_{\boldsymbol{\mathcal{D}}^{(2)}, \theta \circ \psi \circ \varphi} 
\circ 
\left(
\left(
h^{(2)\flat}_{\psi}
\circ
g^{(2)}
\right)_{\varphi}^{\flat}
\circ 
f^{(2)}
\right)^{\flat}
$$

The following chain of equalities holds
\begin{flushleft}
$
\mathrm{pr}^{\llbracket\cdot\rrbracket}_{\boldsymbol{\mathcal{D}}^{(2)}, \theta \circ \psi \circ \varphi} 
\circ 
\left(
h^{(2)\flat}_{\psi \circ \varphi}
\circ
\left(
g^{(2)\flat}_{\varphi}
\circ
f^{(2)}
\right)
\right)^{\flat}
$
\begin{align*}
&=
h^{\llbracket2\rrbracket @}_{\psi \circ \varphi}
\circ
\mathrm{pr}^{\llbracket\cdot\rrbracket}_{\boldsymbol{\mathcal{C}}^{(2)}, \psi \circ \varphi} 
\circ
\left(
g^{(2)\flat}_{\varphi}
\circ
f^{(2)}
\right)^{\flat}
\tag{1}
\\
&=
h^{\llbracket2\rrbracket @}_{\psi \circ \varphi}
\circ
\left(
\mathrm{pr}^{\llbracket\cdot\rrbracket}_{\boldsymbol{\mathcal{C}}^{(2)}, \psi} 
\circ
g^{(2)\flat}
\right)_{\varphi}
\circ
f^{(2)\flat}
\tag{2}
\\
&=
\left(
h^{\llbracket2\rrbracket @}_{\psi}
\right)_{\varphi}
\circ
\left(
\mathrm{pr}^{\llbracket\cdot\rrbracket}_{\boldsymbol{\mathcal{C}}^{(2)}, \psi} 
\circ
g^{(2)\flat}
\right)_{\varphi}
\circ
f^{(2)\flat}
\tag{3}
\\
&=
\left(
h^{\llbracket2\rrbracket @}_{\psi}
\circ
\mathrm{pr}^{\llbracket\cdot\rrbracket}_{\boldsymbol{\mathcal{C}}^{(2)}, \psi} 
\circ
g^{(2)\flat}
\right)_{\varphi}
\circ
f^{(2)\flat}
\tag{4}
\\
&=
\left(
\mathrm{pr}^{\llbracket\cdot\rrbracket}_{\boldsymbol{\mathcal{D}}^{(2)}, \theta \circ \psi}
\circ
\left(
h^{(2)\flat}_{\psi}
\circ
g^{(2)}
\right)^{\flat}
\right)_{\varphi}
\circ
f^{(2)\flat}
\tag{5}
\\
&=
\mathrm{pr}^{\llbracket\cdot\rrbracket}_{\boldsymbol{\mathcal{D}}^{(2)}, \theta \circ \psi \circ \varphi} 
\circ 
\left(
\left(
h^{(2)\flat}_{\psi}
\circ
g^{(2)}
\right)_{\varphi}^{\flat}
\circ 
f^{(2)}
\right)^{\flat}
\tag{6}
\end{align*}
\end{flushleft}

The first equality follows from Proposition~\ref{CDQPthExtComp}, item~(3);
the second equality follows from Proposition~\ref{CDQPthExtComp}, item~(2);
the third equality follows from the fact that, according to Proposition~\ref{PMSetFunc}, $\mathrm{MSet}$ is a contravariant functor;
the fourth equality follows from the fact that, according to Proposition~\ref{PDeltaPhiFunc}, $\Delta_{\varphi}$ is a covariant functor;
the fifth equality follows from Proposition~\ref{CDQPthExtComp}, item~(3);
finally, the last equality follows from Proposition~\ref{CDQPthExtComp}, item~(2).

Thus, Equation (Assoc) follows.

This proves that second-order composition of classes of second-order morphisms between is associative.

This proves that $\mathsf{Rws}^{\llbracket 2 \rrbracket}_{\mathfrak{d}}$ is a category.
\end{proof}

\section{The category $\mathsf{Tw}_{\mathfrak{d}}^{\llbracket 2 \rrbracket}$}

The aim of this section is to define the category of second-order towers and show that it is equivalent to the category of rewriting systems. To begin with, we show that the quotient constructions of second-order equivalent morphisms are equal.

\begin{proposition}\label{PDMorEqv}
Let $\mathbf{f}^{(2)}$ and $\mathbf{f}'^{(2)}$ be second-order morphisms from $\boldsymbol{\mathcal{A}}^{(2)}$ to $\boldsymbol{\mathcal{B}}^{(2)}$. If $\mathbf{f}^{(2)}$ and $\mathbf{f}'^{(2)}$ are second-order equivalent, then 
$\llbracket\mathbf{Pth}_{\boldsymbol{\mathcal{B}}^{(2)}}^{\mathbf{f}^{(2)}}\rrbracket = \llbracket\mathbf{Pth}_{\boldsymbol{\mathcal{B}}^{(2)}}^{\mathbf{f}'^{(2)}}\rrbracket$ and $f^{\llbracket2\rrbracket @} = f'^{\llbracket2\rrbracket @}$.
\end{proposition}

\begin{proof}
Let $\mathbf{f}^{(2)} = (\mathbf{f}^{(1)}, f^{(2)})$ and $\mathbf{f}'^{(2)} = (\mathbf{f}'^{(1)}, f'^{(2)})$ be second-order morphisms from $\boldsymbol{\mathcal{A}}^{(2)}$ to $\boldsymbol{\mathcal{B}}^{(2)}$ such that $\mathbf{f}^{(2)} \cong^{(2)} \mathbf{f}'^{(2)}$. Thus, following Definition~\ref{DDMorEqv}, $\mathbf{f}^{(1)} \cong^{(1)} \mathbf{f}'^{(1)}$ and $\mathrm{pr}^{\llbracket\cdot\rrbracket}_{\boldsymbol{\mathcal{B}}^{(2)}, \varphi} \circ f^{(2)\flat} = \mathrm{pr}^{\llbracket\cdot\rrbracket}_{\boldsymbol{\mathcal{B}}^{(2)}, \varphi} \circ f'^{(2)\flat}$. Therefore, following Definition~\ref{DDQPthExt}, 
$$
f^{\llbracket 2 \rrbracket @}
=
(\mathrm{pr}^{\llbracket\cdot\rrbracket}_{\boldsymbol{\mathcal{B}}^{(2)}, \varphi} \circ f^{(2)\flat})^{\natural}
=
(\mathrm{pr}^{\llbracket\cdot\rrbracket}_{\boldsymbol{\mathcal{B}}^{(2)}, \varphi} \circ f'^{(2)\flat})^{\natural}
=
f'^{\llbracket 2 \rrbracket @}.
$$

Finally, we prove that the many-sorted partial $\Sigma^{\boldsymbol{\mathcal{A}}^{(2)}}$-algebras $\llbracket\mathbf{Pth}_{\boldsymbol{\mathcal{B}}^{(2)}}^{\mathbf{f}^{(2)}}\rrbracket$ and $\llbracket\mathbf{Pth}_{\boldsymbol{\mathcal{B}}^{(2)}}^{\mathbf{f}'^{(2)}}\rrbracket$ are equal. To this end, we consider the different items stated in Proposition~\ref{PDQPthBDCatAlg}.

\textsf{(1)}
The underlying $S$-sorted set of both algebras is $\llbracket\mathrm{Pth}_{\boldsymbol{\mathcal{B}}^{(2)}}\rrbracket_{\varphi}$.

\textsf{(2)}
For every $(\mathbf{s}, s)\in S^{\ast}\times S$ and every operation symbol $\sigma\in\Sigma_{\mathbf{s},s}$, it follows that
$$
\sigma^{\llbracket\mathbf{Pth}_{\boldsymbol{\mathcal{B}}^{(2)}}^{\mathbf{f}^{(2)}}\rrbracket}
=
\sigma^{\mathbf{c}_{\mathfrak{d}}^{\ast}(\llbracket\mathbf{Pth}_{\boldsymbol{\mathcal{B}}^{(2)}}^{(0,2)}\rrbracket)}
=
\sigma^{\llbracket\mathbf{Pth}_{\boldsymbol{\mathcal{B}}^{(2)}}^{\mathbf{f}'^{(2)}}\rrbracket}.
$$

\textsf{(3)}
For every sort $s$ in $S$ and every rewrite rule $\mathfrak{p} \in \mathcal{A}_{s}^{(1)}$, the following chain of equalities holds
\allowdisplaybreaks
\begin{align*}
\mathfrak{p}^{\llbracket\mathbf{Pth}_{\boldsymbol{\mathcal{B}}^{(2)}}^{\mathbf{f}^{(2)}}\rrbracket}
&=
\mathrm{pr}^{\llbracket\cdot\rrbracket}_{\boldsymbol{\mathcal{B}}^{(2)}, \varphi(s)}
\circ
f^{(2)\flat}_{(s)}
\circ
\mathrm{ech}^{(2, \mathcal{A}^{(1)})}_{\boldsymbol{\mathcal{A}}^{(2)}, s}
\left(
\mathfrak{p}^{(1)}
\right)
\tag{1}
\\
&=
\mathrm{pr}^{\llbracket\cdot\rrbracket}_{\boldsymbol{\mathcal{B}}^{(2)}, \varphi(s)}
\circ
f'^{(2)\flat}_{(s)}
\circ
\mathrm{ech}^{(2, \mathcal{A}^{(1)})}_{\boldsymbol{\mathcal{A}}^{(2)}, s}
\left(
\mathfrak{p}^{(1)}
\right)
\tag{2}
\\
&=
\mathfrak{p}^{\llbracket\mathbf{Pth}_{\boldsymbol{\mathcal{B}}^{(2)}}^{\mathbf{f}'^{(2)}}\rrbracket}.
\tag{3}
\end{align*}

The first equality unravels the definition of the interpretation of the constant symbol $\mathfrak{p}$ in the partial $\Sigma^{\boldsymbol{\mathcal{A}}^{(2)}}$-algebra $\llbracket\mathbf{Pth}_{\boldsymbol{\mathcal{B}}^{(2)}}^{\mathbf{f}^{(2)}}\rrbracket$;
the second equality follows from the fact that $\mathbf{f}^{(2)} \cong^{(2)} \mathbf{f}'^{(2)}$;
finally, the last equality recovers the definition of the interpretation of the constant symbol $\mathfrak{p}$ in the partial $\Sigma^{\boldsymbol{\mathcal{A}}^{(2)}}$-algebra $\llbracket\mathbf{Pth}_{\boldsymbol{\mathcal{B}}^{(2)}}^{\mathbf{f}'^{(2)}}\rrbracket$;

\textsf{(4)}
For every sort $s$ in $S$, it follows that
\begin{align*}
\mathrm{sc}_{s}^{0\llbracket\mathbf{Pth}_{\boldsymbol{\mathcal{B}}^{(2)}}^{\mathbf{f}^{(2)}}\rrbracket}
&=
\mathrm{sc}_{\varphi(s)}^{0\llbracket\mathbf{Pth}_{\boldsymbol{\mathcal{B}}^{(2)}}\rrbracket}
=
\mathrm{sc}_{s}^{0\llbracket\mathbf{Pth}_{\boldsymbol{\mathcal{B}}^{(2)}}^{\mathbf{f}'^{(2)}}\rrbracket};
\\
\mathrm{tg}_{s}^{0\llbracket\mathbf{Pth}_{\boldsymbol{\mathcal{B}}^{(2)}}^{\mathbf{f}^{(2)}}\rrbracket}
&=
\mathrm{tg}_{\varphi(s)}^{0\llbracket\mathbf{Pth}_{\boldsymbol{\mathcal{B}}^{(2)}}\rrbracket}
=
\mathrm{tg}_{s}^{0\llbracket\mathbf{Pth}_{\boldsymbol{\mathcal{B}}^{(2)}}^{\mathbf{f}'^{(2)}}\rrbracket}.
\end{align*}

\textsf{(5)}
Similarly, for every sort $s$ in $S$, it follows that
$$
\circ_{s}^{0\llbracket\mathbf{Pth}_{\boldsymbol{\mathcal{B}}^{(2)}}^{\mathbf{f}^{(2)}}\rrbracket}
=
\circ_{\varphi(s)}^{0\llbracket\mathbf{Pth}_{\boldsymbol{\mathcal{B}}^{(2)}}\rrbracket}
=
\circ_{s}^{0\llbracket\mathbf{Pth}_{\boldsymbol{\mathcal{B}}^{(2)}}^{\mathbf{f}'^{(2)}}\rrbracket}.
$$

\textsf{(6)}
For every sort $s$ in $S$ and every rewrite rule $\mathfrak{p}^{(2)} \in \mathcal{A}_{s}^{(2)}$, the following chain of equalities holds
\allowdisplaybreaks
\begin{align*}
\mathfrak{p}^{(2)\llbracket\mathbf{Pth}_{\boldsymbol{\mathcal{B}}^{(2)}}^{\mathbf{f}^{(2)}}\rrbracket}
&=
\mathrm{pr}^{\llbracket\cdot\rrbracket}_{\boldsymbol{\mathcal{B}}^{(2)}, \varphi(s)}
\circ
f^{(2)\flat}_{(s)}
\left(
\mathfrak{p}^{(2)}
\right)
\tag{1}
\\
&=
\mathrm{pr}^{\llbracket\cdot\rrbracket}_{\boldsymbol{\mathcal{B}}^{(2)}, \varphi(s)}
\circ
f^{(2)\flat}_{(s)}
\circ
\mathrm{ech}^{(2, \mathcal{A}^{(2)})}_{\boldsymbol{\mathcal{A}}^{(2)}, s}
\left(
\mathfrak{p}^{(2)}
\right)
\tag{2}
\\
&=
\mathrm{pr}^{\llbracket\cdot\rrbracket}_{\boldsymbol{\mathcal{B}}^{(2)}, \varphi(s)}
\circ
f'^{(2)\flat}_{(s)}
\circ
\mathrm{ech}^{(2, \mathcal{A}^{(2)})}_{\boldsymbol{\mathcal{A}}^{(2)}, s}
\left(
\mathfrak{p}^{(2)}
\right)
\tag{3}
\\
&=
\mathrm{pr}^{\llbracket\cdot\rrbracket}_{\boldsymbol{\mathcal{B}}^{(2)}, \varphi(s)}
\circ
f'^{(2)\flat}_{(s)}
\left(
\mathfrak{p}^{(2)}
\right)
\tag{4}
\\
&=
\mathfrak{p}^{(2)\llbracket\mathbf{Pth}_{\boldsymbol{\mathcal{B}}^{(2)}}^{\mathbf{f}^{(2)}}\rrbracket}.
\tag{5}
\end{align*}

The first equality unravels the definition of the interpretation of the constant symbol $\mathfrak{p}^{(2)}$ in the partial $\Sigma^{\boldsymbol{\mathcal{A}}^{(2)}}$-algebra $\llbracket\mathbf{Pth}_{\boldsymbol{\mathcal{B}}^{(2)}}^{\mathbf{f}^{(2)}}\rrbracket$;
the second equality follows from the definition of the path extension mapping, introduced in Proposition~\ref{PPthExt};
the third equality follows from the fact that $\mathbf{f}^{(2)} \cong^{(2)} \mathbf{f}'^{(2)}$;
the fourth equality follows from the definition of the path extension mapping, introduced in Proposition~\ref{PPthExt};
finally, the last equality recovers the definition of the interpretation of the constant symbol $\mathfrak{p}$ in the partial $\Sigma^{\boldsymbol{\mathcal{A}}^{(2)}}$-algebra $\llbracket\mathbf{Pth}_{\boldsymbol{\mathcal{B}}^{(2)}}^{\mathbf{f}'^{(2)}}\rrbracket$;

\textsf{(7)}
For every sort $s$ in $S$, it follows that
\begin{align*}
\mathrm{sc}_{s}^{1\llbracket\mathbf{Pth}_{\boldsymbol{\mathcal{B}}^{(2)}}^{\mathbf{f}^{(2)}}\rrbracket}
&=
\mathrm{sc}_{\varphi(s)}^{1\llbracket\mathbf{Pth}_{\boldsymbol{\mathcal{B}}^{(2)}}\rrbracket}
=
\mathrm{sc}_{s}^{1\llbracket\mathbf{Pth}_{\boldsymbol{\mathcal{B}}^{(2)}}^{\mathbf{f}'^{(2)}}\rrbracket};
\\
\mathrm{tg}_{s}^{1\llbracket\mathbf{Pth}_{\boldsymbol{\mathcal{B}}^{(2)}}^{\mathbf{f}^{(2)}}\rrbracket}
&=
\mathrm{tg}_{\varphi(s)}^{1\llbracket\mathbf{Pth}_{\boldsymbol{\mathcal{B}}^{(2)}}\rrbracket}
=
\mathrm{tg}_{s}^{1\llbracket\mathbf{Pth}_{\boldsymbol{\mathcal{B}}^{(2)}}^{\mathbf{f}'^{(2)}}\rrbracket}.
\end{align*}

\textsf{(8)}
Similarly, for every sort $s$ in $S$, it follows that
$$
\circ_{s}^{1\llbracket\mathbf{Pth}_{\boldsymbol{\mathcal{B}}^{(2)}}^{\mathbf{f}^{(2)}}\rrbracket}
=
\circ_{\varphi(s)}^{1\llbracket\mathbf{Pth}_{\boldsymbol{\mathcal{B}}^{(2)}}\rrbracket}
=
\circ_{s}^{1\llbracket\mathbf{Pth}_{\boldsymbol{\mathcal{B}}^{(2)}}^{\mathbf{f}'^{(2)}}\rrbracket}.
$$
\end{proof}

We now define, given a second-order rewriting system, its associated second-order tower and the notion of morphism of second-order towers.

\begin{definition}\label{DTw2}
Given a second-order many-sorted rewriting system $\boldsymbol{\mathcal{A}}^{(2)} = (\boldsymbol{\mathcal{A}}^{(1)}, \mathcal{A}^{(2)})$, the \emph{second-order tower} associated to $\boldsymbol{\mathcal{A}}^{(2)}$ is 
$$
\mathbb{A}^{(2)}
=
\left(
\boldsymbol{\mathcal{A}}^{(2)}, 
\mathbf{T}_{\boldsymbol{\mathcal{E}}^{\boldsymbol{\mathcal{A}}^{(2)}}}\left(\mathbf{Pth}_{\boldsymbol{\mathcal{A}}^{(2)}}\right)
\right)
$$
where $\mathbf{T}_{\boldsymbol{\mathcal{E}}^{\boldsymbol{\mathcal{A}}^{(2)}}}(\mathbf{Pth}_{\boldsymbol{\mathcal{A}}^{(2)}})$ is the free $\Sigma^{\boldsymbol{\mathcal{A}}^{(2)}}$-algebra on the QE-variety $\mathcal{V}(\boldsymbol{\mathcal{E}}^{\boldsymbol{\mathcal{A}}^{(2)}})$. 

A \emph{morphisms of second-order towers}from $\mathbb{A}^{(2)}=(\boldsymbol{\mathcal{A}}^{(2)}, \mathbf{T}_{\boldsymbol{\mathcal{E}}^{\boldsymbol{\mathcal{A}}^{(2)}}}(\mathbf{Pth}_{\boldsymbol{\mathcal{A}}^{(2)}})$ to $\mathbb{B}^{(2)}=(\boldsymbol{\mathcal{B}}^{(2)}, \mathbf{T}_{\boldsymbol{\mathcal{E}}^{\boldsymbol{\mathcal{B}}^{(2)}}}(\mathbf{Pth}_{\boldsymbol{\mathcal{B}}^{(2)}})$, or simply a \emph{second-order tower morphism}, is an ordered triple $(\mathbb{A}^{(2)}, \mathbf{f}^{\llbracket2\rrbracket @}, \mathbb{B}^{(2)})$, denoted by $\mathbf{f}^{\llbracket2\rrbracket @} \colon \mathbb{A}^{(2)} \mor \mathbb{B}^{(2)}$ for short, in which $\mathbf{f}^{\llbracket2\rrbracket @}$ is the ordered pair $(\llbracket\mathbf{f}^{(2)}\rrbracket, f^{\llbracket2\rrbracket @})$ where $\llbracket\mathbf{f}^{(2)}\rrbracket$ is the second-order morphism class of a second-order morphism $\mathbf{f}^{(2)}$ from $\boldsymbol{\mathcal{A}}^{(2)}$ to $\boldsymbol{\mathcal{B}}^{(2)}$ and $f^{\llbracket2\rrbracket @}$ is the $\Sigma^{\boldsymbol{\mathcal{A}}^{(2)}}$-homomorphism from $\mathbf{T}_{\boldsymbol{\mathcal{E}}^{\boldsymbol{\mathcal{A}}^{(2)}}}(\mathbf{Pth}_{\boldsymbol{\mathcal{A}}^{(2)}})$ to $\mathbf{T}_{\boldsymbol{\mathcal{E}}^{\boldsymbol{\mathcal{B}}^{(2)}}}^{\mathbf{f}^{(2)}}(\mathbf{Pth}_{\boldsymbol{\mathcal{B}}^{(2)}})$ introduced in Definition~\ref{DDQPthExt}.

Note that the mapping introduced in Definition~\ref{DDQPthExt}, is a $\Sigma^{\boldsymbol{\mathcal{A}}^{(2)}}$-homomorphism from $\llbracket \mathbf{Pth}_{\boldsymbol{\mathcal{A}}^{(2)}}\rrbracket$ to $\llbracket \mathbf{Pth}_{\boldsymbol{\mathcal{B}}^{(2)}}^{\mathbf{f}^{(2)}}\rrbracket$. However, taking into account Theorems~\ref{TDPthFree} and \ref{TDPthBFreeB}, we obtain a unique $\Sigma^{\boldsymbol{\mathcal{A}}^{(2)}}$-homomorphism, also denoted by $f^{\llbracket2\rrbracket @}$ along this subsection, from $\mathbf{T}_{\boldsymbol{\mathcal{E}}^{\boldsymbol{\mathcal{A}}^{(2)}}}(\mathbf{Pth}_{\boldsymbol{\mathcal{A}}^{(2)}})$ to $\mathbf{T}_{\boldsymbol{\mathcal{E}}^{\boldsymbol{\mathcal{B}}^{(2)}}}^{\mathbf{f}^{(2)}}(\mathbf{Pth}_{\boldsymbol{\mathcal{B}}^{(2)}})$.
\end{definition}

We define below the notions of second-order identity morphism at a second-order tower and second-order composite morphism of second-order towers.

\begin{definition}
\label{DIdTw2}
Let $\mathbb{A}^{(2)}=(\boldsymbol{\mathcal{A}}^{(2)}, \mathbf{T}_{\boldsymbol{\mathcal{E}}^{\boldsymbol{\mathcal{A}}^{(2)}}}(\mathbf{Pth}_{\boldsymbol{\mathcal{A}}^{(2)}}))$ be the second-order tower associated to the second-order many-sorted rewriting system $\boldsymbol{\mathcal{A}}^{(2)}$. The ordered pair $(\llbracket\mathrm{id}^{\boldsymbol{\mathcal{A}}^{(2)}}\rrbracket, \mathrm{id}^{\mathbf{T}_{\boldsymbol{\mathcal{E}}^{\boldsymbol{\mathcal{A}}^{(2)}}}(\mathbf{Pth}_{\boldsymbol{\mathcal{A}}^{(2)}})})$, denoted by $\mathrm{id}^{\mathbb{A}^{(2)}}$, is the second-order identity morphism at $\mathbb{A}^{(2)}$. Recall that, $\llbracket\mathrm{id}^{\boldsymbol{\mathcal{A}}^{(2)}}\rrbracket$ is the second-order morphism class of the second-order identity morphism at $\boldsymbol{\mathcal{A}}^{(2)}$ defined as
$$
\mathrm{id}^{\boldsymbol{\mathcal{A}}^{(2)}}
=
\left(
\mathrm{id}^{\boldsymbol{\mathcal{A}}^{(1)}},
\mathrm{ech}^{(2,\mathcal{A}^{(2)})}_{\boldsymbol{\mathcal{A}}^{(2)}}
\right)
$$
and, according to Proposition~\ref{PDQPthExtEch}, $\mathrm{id}^{\mathbf{T}_{\boldsymbol{\mathcal{E}}^{\boldsymbol{\mathcal{A}}^{(2)}}}(\mathbf{Pth}_{\boldsymbol{\mathcal{A}}^{(2)}})}$ is its quotient path extension mapping.
\end{definition}

\begin{definition}
\label{DCompTw2}
Let $\mathbf{f}^{\llbracket2\rrbracket @}=(\llbracket\mathbf{f}^{(2)}\rrbracket, f^{\llbracket2\rrbracket @})$ be a second-order tower morphism from $\mathbb{A}^{(2)}$ to $\mathbb{B}^{(2)}$ and $\mathbf{g}^{\llbracket2\rrbracket @}=(\llbracket\mathbf{g}^{(2)}\rrbracket, g^{\llbracket2\rrbracket @})$ a second-order tower morphism from $\mathbb{B}^{(2)}$ to $\mathbb{C}^{(2)}$ where $\mathbb{A}^{(2)}$, $\mathbb{B}^{(2)}$ and $\mathbb{C}^{(2)}$ stand for
$$
\mathbb{A}^{(2)}=(\boldsymbol{\mathcal{A}}^{(2)}, \mathbf{T}_{\boldsymbol{\mathcal{E}}^{\boldsymbol{\mathcal{A}}^{(2)}}}(\mathbf{Pth}_{\boldsymbol{\mathcal{A}}^{(2)}})), \,
\mathbb{B}^{(2)}=(\boldsymbol{\mathcal{B}}^{(2)}, \mathbf{T}_{\boldsymbol{\mathcal{E}}^{\boldsymbol{\mathcal{B}}^{(2)}}}(\mathbf{Pth}_{\boldsymbol{\mathcal{B}}^{(2)}}))
$$
$$
\mbox{ and }
\mathbb{C}^{(2)}=(\boldsymbol{\mathcal{C}}^{(2)}, \mathbf{T}_{\boldsymbol{\mathcal{E}}^{\boldsymbol{\mathcal{C}}^{(2)}}}(\mathbf{Pth}_{\boldsymbol{\mathcal{C}}^{(2)}})),
\mbox{ respectively.}
$$
The \emph{second-order composite morphism} $\mathbf{g}^{\llbracket2\rrbracket @}\circ\mathbf{f}^{\llbracket2\rrbracket @}$, from $\mathbb{A}^{(2)}$ to $\mathbb{C}^{(2)}$, is
$$
\mathbf{g}^{\llbracket2\rrbracket @}\circ\mathbf{f}^{\llbracket2\rrbracket @}=
\left(
\left\llbracket\mathbf{g}^{(2)}\right\rrbracket\circ\left\llbracket\mathbf{f}^{(2)}\right\rrbracket, 
g^{\llbracket2\rrbracket @}_{\varphi} \circ f^{\llbracket2\rrbracket @}
\right).
$$
Recall that the composition $\llbracket\mathbf{g}^{(2)}\rrbracket \circ \llbracket\mathbf{f}^{(2)}\rrbracket$ is defined to be the second-order morphism class $\llbracket\mathbf{g}^{(2)} \circ \mathbf{f}^{(2)}\rrbracket$ and, according to Corollary~\ref{CDQPthExtComp}, $g^{\llbracket2\rrbracket @}_{\varphi} \circ f^{\llbracket2\rrbracket @}$ is its quotient second-order path extension mapping.
\end{definition}

We show that second-order towers and second-order morphisms between them constitute a category.

\begin{proposition}\label{PTw2Cat}
The second-order towers together with the second-order morphisms between second-order towers constitute a category, denoted by $\mathsf{Tw}^{(2)}_{\mathfrak{d}}$.
\end{proposition}

\begin{proof}
That domains and codomains respect identities and compositions follows from the definitions of second-order identity morphism and second-order composition of second-order morphisms introduced in Definition~\ref{DCompTw2}. Thus, all that remains to be proven is that the second-order identity morphism at a second-order tower acts as a unit element and that the second-order composition of morphisms of second-order towers is associative.

{\sffamily Unit element.}


Let $\mathbf{f}^{\llbracket2\rrbracket @}=(\llbracket\mathbf{f}^{(2)}\rrbracket, f^{\llbracket2\rrbracket @})$ be a second-order tower morphism from $\mathbb{A}^{(2)}$ to $\mathbb{B}^{(2)}$. We need to prove that 
\begin{align*}
\mathbf{f}^{\llbracket2\rrbracket @} \circ \mathrm{id}^{\mathbb{A}^{(2)}} &= \mathbf{f}^{\llbracket2\rrbracket @}
&&\mbox{ and }&
\mathrm{id}^{\mathbb{B}^{(2)}} \circ \mathbf{f}^{\llbracket2\rrbracket @} &= \mathbf{f}^{\llbracket2\rrbracket @}.
\tag{Id}
\end{align*}

We will only prove the left hand side of (Id). The right hand side of (Id) is done similarly. Let us prove that $\mathbf{f}^{\llbracket2\rrbracket @} \circ \mathrm{id}^{\mathbb{A}^{(2)}} = \mathbf{f}^{\llbracket2\rrbracket @}$.

The following chain of equalities holds
\allowdisplaybreaks
\begin{align*}
\mathbf{f}^{\llbracket2\rrbracket @} \circ \mathrm{id}^{\mathbb{A}^{(2)}}
&=
\left( \left\llbracket\mathbf{f}^{(2)}\right\rrbracket, f^{\llbracket2\rrbracket @} \right) \circ \left( \left\llbracket\mathrm{id}^{\boldsymbol{\mathcal{A}}^{(2)}}\right\rrbracket, \mathrm{id}^{\mathbf{T}_{\boldsymbol{\mathcal{E}}^{\boldsymbol{\mathcal{A}}^{(2)}}}(\mathbf{Pth}_{\boldsymbol{\mathcal{A}}^{(2)}})} \right)
\tag{1}
\\
&=
\left(
\left\llbracket\mathbf{f}^{(2)}\right\rrbracket \circ \left\llbracket\mathrm{id}^{\boldsymbol{\mathcal{A}}^{(2)}}\right\rrbracket
,
f^{\llbracket2\rrbracket @}_{\mathrm{id}^{S}} \circ \mathrm{id}^{\mathbf{T}_{\boldsymbol{\mathcal{E}}^{\boldsymbol{\mathcal{A}}^{(2)}}}(\mathbf{Pth}_{\boldsymbol{\mathcal{A}}^{(2)}})}
\right)
\tag{2}
\\
&=
\left(
\left\llbracket\mathbf{f}^{(2)}\right\rrbracket
,
f^{\llbracket2\rrbracket @} \circ \mathrm{id}^{\mathbf{T}_{\boldsymbol{\mathcal{E}}^{\boldsymbol{\mathcal{A}}^{(2)}}}(\mathbf{Pth}_{\boldsymbol{\mathcal{A}}^{(2)}})}
\right)
\tag{3}
\\
&=
\left(
\left\llbracket\mathbf{f}^{(2)}\right\rrbracket
,
f^{\llbracket2\rrbracket @}
\right)
\tag{4}
\\
&=
\mathbf{f}^{\llbracket2\rrbracket@}
\tag{5}
\end{align*}

The first equality unravels the definitions of the second-order tower morphisms $\mathbf{f}^{(2)@}$ and $\mathrm{id}^{\mathbb{A}^{(2)}}$;
the second equality unravels the definition of second-order composition of second-order tower morphisms, introduced in Definition~\ref{DCompTw2};
the third equality follows from the fact that, according to Proposition~\ref{PMSetFunc}, $\mathrm{MSet}$ is a contravariant functor;
the fourth equality follows from the fact that $f^{\llbracket2\rrbracket @} \circ \mathrm{id}^{\mathbf{T}_{\boldsymbol{\mathcal{E}}^{\boldsymbol{\mathcal{A}}^{(2)}}}(\mathbf{Pth}_{\boldsymbol{\mathcal{A}}^{(2)}})}$ is equal to $f^{\llbracket2\rrbracket @}$;
finally, the last equality recovers the definition of $\mathbf{f}^{\llbracket2\rrbracket @}$.

Thus, Equation (Id) follows.

This proves that second-order identities at second-order towers are the unit element for second-order composition of morphisms.

{\sffamily Associativity.}

Let $\mathbb{A}^{(2)}$, $\mathbb{B}^{(2)}$, $\mathbb{C}^{(2)}$ and $\mathbb{D}^{(2)}$ be four second-order towers and let $\mathbf{f}^{\llbracket2\rrbracket @}=(\llbracket\mathbf{f}^{(2)}\rrbracket, f^{\llbracket2\rrbracket @})$, $\mathbf{g}^{\llbracket2\rrbracket @}=(\llbracket\mathbf{g}^{(2)}\rrbracket, g^{\llbracket2\rrbracket @})$ and $\mathbf{h}^{\llbracket2\rrbracket @}=(\llbracket\mathbf{h}^{(2)}\rrbracket, h^{\llbracket2\rrbracket @})$ be second-order tower morphisms of the form
$$
\mathbf{f}^{\llbracket2\rrbracket @} \colon \mathbb{A}^{(2)} \mor \mathbb{B}^{(2)}, \,
\mathbf{g}^{\llbracket2\rrbracket @} \colon \mathbb{B}^{(2)} \mor \mathbb{C}^{(2)}
\mbox{ and }
\mathbf{h}^{\llbracket2\rrbracket @} \colon \mathbb{C}^{(2)} \mor \mathbb{D}^{(2)}
$$
where $\mathbf{f}^{(2)}$, $\mathbf{g}^{(2)}$ and $\mathbf{h}^{(2)}$ stand for
$$
\mathbf{f}^{(2)}=(\varphi, c, (f^{(i)})_{i\in 3}), \,
\mathbf{g}^{(2)}=(\psi, d, (g^{(i)})_{i\in 3})
\mbox{ and }
\mathbf{h}^{(2)}=(\theta, e, (h^{(i)})_{i\in 3}),
\mbox{ respectively.}
$$
We need to prove that
\begin{equation}
\mathbf{h}^{\llbracket2\rrbracket @} \circ \left(\mathbf{g}^{\llbracket2\rrbracket @} \circ \mathbf{f}^{\llbracket2\rrbracket @}\right)
=
\left(\mathbf{h}^{\llbracket2\rrbracket @} \circ \mathbf{g}^{\llbracket2\rrbracket @}\right) \circ \mathbf{f}^{\llbracket2\rrbracket @}.
\tag{Assoc}
\end{equation}

The following chain of equalities holds
\allowdisplaybreaks
\begin{align*}
&
\mathbf{h}^{\llbracket2\rrbracket @} \circ \left(\mathbf{g}^{\llbracket2\rrbracket @} \circ \mathbf{f}^{\llbracket2\rrbracket @}\right)
\\
&=
\left( \llbracket\mathbf{h}^{(2)}\rrbracket,h^{\llbracket2\rrbracket @} \right) \circ \left(\left( \llbracket\mathbf{g}^{(2)}\rrbracket, g^{\llbracket2\rrbracket @} \right) \circ \left( \llbracket\mathbf{f}^{(2)}\rrbracket, f^{\llbracket2\rrbracket @} \right) \right)
\tag{1}
\\
&=
\left(
\llbracket\mathbf{h}^{(2)}\rrbracket \circ \left(\llbracket\mathbf{g}^{(2)}\rrbracket \circ \llbracket\mathbf{f}^{(2)}\rrbracket\right)
,
h^{\llbracket2\rrbracket @}_{\psi \circ \varphi} \circ g^{\llbracket2\rrbracket @}_{\varphi} \circ f^{\llbracket2\rrbracket @}
\right)
\tag{2}
\\
&=
\left(
\llbracket\mathbf{h}^{(2)}\rrbracket \circ \left(\llbracket\mathbf{g}^{(2)}\rrbracket \circ \llbracket\mathbf{f}^{(2)}\rrbracket\right)
,
(h^{\llbracket2\rrbracket @}_{\psi})_{\varphi} \circ g^{\llbracket2\rrbracket @}_{\varphi} \circ f^{\llbracket2\rrbracket @}
\right)
\tag{3}
\\
&=
\left(
\left(\llbracket\mathbf{h}^{(2)}\rrbracket \circ \llbracket\mathbf{g}^{(2)}\rrbracket\right) \circ \llbracket\mathbf{f}^{(2)}\rrbracket
,
\left(h^{\llbracket2\rrbracket @}_{\psi} \circ g^{\llbracket2\rrbracket @}\right)_{\varphi} \circ f^{\llbracket2\rrbracket @}
\right)
\tag{4}
\\
&=
\left( \left( \llbracket\mathbf{h}^{(2)}\rrbracket,h^{\llbracket2\rrbracket @} \right) \circ \left( \llbracket\mathbf{g}^{(2)}\rrbracket, g^{\llbracket2\rrbracket @} \right) \right) \circ \left( \llbracket\mathbf{f}^{(2)}\rrbracket, f^{\llbracket2\rrbracket @} \right)
\tag{5}
\\
&=
\left(\mathbf{h}^{\llbracket2\rrbracket @} \circ \mathbf{g}^{\llbracket2\rrbracket @}\right) \circ \mathbf{f}^{\llbracket2\rrbracket @}.
\tag{6}
\end{align*}

The first equality unravels the definitions of the second-order tower morphisms $\mathbf{f}^{\llbracket2\rrbracket @}$, $\mathbf{g}^{\llbracket2\rrbracket @}$ and $\mathbf{h}^{\llbracket2\rrbracket @}$;
the second equality unravels the definition of second-order composition of second-order tower morphisms, introduced in Definition~\ref{DCompTw2};
the third equality follows from the fact that, by Proposition~\ref{PMSetFunc}, $\mathrm{MSet}$ is a contravariant functor;
the fourth equality follows from the fact that, by Proposition~\ref{PDeltaPhiFunc}, $\Delta_{\varphi}$ is a covariant functor;
the fifth equality recovers the definition of second-order composition of second-order tower morphisms, introduced in Definition~\ref{DCompTw2};
finally, the last equality recovers the definitions of the second-order tower morphisms $\mathbf{f}^{\llbracket2\rrbracket @}$, $\mathbf{g}^{\llbracket2\rrbracket @}$ and $\mathbf{h}^{\llbracket2\rrbracket @}$.

Thus, Equation (Assoc) follows.

This proves that second-order composition of second-order tower morphisms is associative.

This shows that $\mathsf{Tw}_{\mathfrak{d}}^{\llbracket 2 \rrbracket}$ is a category.
\end{proof}

We will prove bellow that the categories of second-order many-sorted rewriting systems and second-order towers are isomorphic. In particular, we define the assignments $V^{(2)}$ and $U^{(2)}$, we prove that they are functors between the respective categories and that they are mutually inverse functors.

\begin{definition}\label{DV2}
We let $V^{(2)}$ stand for the assignment from $\mathsf{Rws}_{\mathfrak{d}}^{\llbracket 2 \rrbracket}$ to $\mathsf{Tw}_{\mathfrak{d}}^{\llbracket 2 \rrbracket}$ defined as follows:
\begin{enumerate}
\item
for every second-order many-sorted rewriting system $\boldsymbol{\mathcal{A}}^{(2)}$, $V^{(2)}(\boldsymbol{\mathcal{A}}^{(2)})$ is the associated second-order tower $\mathbb{A}^{(2)} = (\boldsymbol{\mathcal{A}}^{(2)}, \mathbf{T}_{\boldsymbol{\mathcal{E}}^{\boldsymbol{\mathcal{A}}^{(2)}}}(\mathbf{Pth}_{\boldsymbol{\mathcal{A}}^{(2)}}))$, and
\item
for every second-order morphism $\mathbf{f}^{(2)} \colon \boldsymbol{\mathcal{A}}^{(2)} \mor \boldsymbol{\mathcal{B}}^{(2)}$, $V^{(2)}(\llbracket\mathbf{f}^{(2)}\rrbracket)$ is the second-order tower morphism $\mathbf{f}^{\llbracket2\rrbracket @}=(\llbracket\mathbf{f}^{(2)}\rrbracket, f^{\llbracket2\rrbracket @}) \colon \mathbb{A}^{(2)} \mor \mathbb{B}^{(2)}$.
\end{enumerate}
\end{definition}

\begin{proposition}\label{PV2Fun}
The assignment $V^{(2)}$ from $\mathsf{Rws}_{\mathfrak{d}}^{\llbracket 2 \rrbracket}$ to $\mathsf{Tw}_{\mathfrak{d}}^{\llbracket 2 \rrbracket}$ is a covariant functor.
\end{proposition}

\begin{proof}
That $V^{(2)}$ maps objects and morphisms of $\mathsf{Rws}_{\mathfrak{d}}^{\llbracket 2 \rrbracket}$ to objects and morphisms of $\mathsf{Tw}_{\mathfrak{d}}^{\llbracket 2 \rrbracket}$ follows from the definition of the assignment. Therefore, all that remains to be proven is that $V^{(2)}$ preserves identities and compositions.

{\sffamily $V^{(2)}$ preserves identities.}

Let $\boldsymbol{\mathcal{A}}^{(2)}$ be an object in $\mathsf{Rws}_{\mathfrak{d}}^{\llbracket 2 \rrbracket}$. We need to prove that
$$
V^{(2)} \left(\left\llbracket\mathrm{id}^{\boldsymbol{\mathcal{A}}^{(2)}}\right\rrbracket\right)
=
\mathrm{id}^{V^{(2)}(\boldsymbol{\mathcal{A}}^{(2)})}.
$$

The following chain of equalities holds
\allowdisplaybreaks
\begin{align*}
V^{(2)} \left(\left\llbracket\mathrm{id}^{\boldsymbol{\mathcal{A}}^{(2)}}\right\rrbracket\right)
&=
\left( \left\llbracket\mathrm{id}^{\boldsymbol{\mathcal{A}}^{(2)}}\right\rrbracket, \mathrm{ech}^{(2,\mathcal{A}^{(2)})@}_{\boldsymbol{\mathcal{A}}^{(2)}}\right)
\tag{1}
\\
&=
\left( \left\llbracket\mathrm{id}^{\boldsymbol{\mathcal{A}}^{(2)}}\right\rrbracket, \mathrm{id}^{\mathbf{T}_{\boldsymbol{\mathcal{E}}^{\boldsymbol{\mathcal{A}}^{(2)}}}(\mathbf{Pth}_{\boldsymbol{\mathcal{A}}^{(2)}})}\right)
\tag{2}
\\
&=
\mathrm{id}^{(\boldsymbol{\mathcal{A}}^{(2)}, \mathbf{T}_{\boldsymbol{\mathcal{E}}^{\boldsymbol{\mathcal{A}}^{(2)}}}(\mathbf{Pth}_{\boldsymbol{\mathcal{A}}^{(2)}}))}
\tag{3}
\\
&=
\mathrm{id}^{V^{(2)}(\boldsymbol{\mathcal{A}}^{2})}.
\tag{4}
\end{align*}
The first equality unravels the definition of the assignment $V^{(2)}$ on morphisms, introduced in Definition~\ref{DV2};
the second equality follows from Proposition~\ref{PDQPthExtEch}; 
the third equality recovers the definition of the second-order identity tower morphism at the second-order tower $(\boldsymbol{\mathcal{A}}^{(2)}, \mathbf{T}_{\boldsymbol{\mathcal{E}}^{\boldsymbol{\mathcal{A}}^{(2)}}}(\mathbf{Pth}_{\boldsymbol{\mathcal{A}}^{(2)}}))$, introduced in Definition~\ref{DIdTw2}; 
finally, the last equality recovers the definition of the assignment $V^{(2)}$ on objects.

This proves that $V^{(2)}$ preserves identities.

{\sffamily $V^{(2)}$ Preserves compositions.}

Let $\mathbf{f}^{(2)}=(\varphi, c, (f^{(i)})_{i\in 3}) \colon \boldsymbol{\mathcal{A}}^{(2)} \mor \boldsymbol{\mathcal{B}}^{(2)}$ and $\mathbf{g}^{(2)}=(\psi, d, (g^{(i)})_{i \in 3}) \colon \boldsymbol{\mathcal{B}}^{(2)} \mor \boldsymbol{\mathcal{C}}^{(2)}$ be second-order morphisms. We need to prove that
$$
V^{(2)} \left(\left\llbracket\mathbf{g}^{(2)}\right\rrbracket
\circ
\left\llbracket\mathbf{f}^{(2)}\right\rrbracket\right)
=
V^{(2)} \left(\left\llbracket\mathbf{g}^{(2)}\right\rrbracket\right) \circ V^{(2)} \left(\left\llbracket\mathbf{g}^{(2)}\right\rrbracket\right).
$$

The following chain of equalities holds
\begin{align*}
V^{(2)} \left(\left\llbracket\mathbf{g}^{(2)}\right\rrbracket
\circ
\left\llbracket\mathbf{f}^{(2)}\right\rrbracket\right)
&=
V^{(2)} \left(\left\llbracket\mathbf{g}^{(2)}\circ\mathbf{f}^{(2)}\right\rrbracket\right)
\tag{1}
\\
&=
\left(
\left\llbracket\mathbf{g}^{(2)}\circ\mathbf{f}^{(2)}\right\rrbracket,
\left(
g^{(2)\flat}_{\varphi} \circ f^{(2)}
\right)^{@}
\right)
\tag{2}
\\
&=
\left(
\left\llbracket\mathbf{g}^{(2)}\right\rrbracket\circ\left\llbracket\mathbf{f}^{(2)}\right\rrbracket,
\left(
g^{(2)\flat}_{\varphi} \circ f^{(2)}
\right)^{@}
\right)
\tag{3}
\\
&=
\left(
\left\llbracket\mathbf{g}^{(2)}\right\rrbracket\circ\left\llbracket\mathbf{f}^{(2)}\right\rrbracket,
g^{\llbracket2\rrbracket @}_{\varphi} \circ f^{\llbracket2\rrbracket @}
\right)
\tag{4}
\\
&=
\mathbf{g}^{\llbracket2\rrbracket @}
\circ
\mathbf{f}^{\llbracket2\rrbracket @}
\tag{5}
\\
&=
V^{(2)} \left(\left\llbracket\mathbf{g}^{(2)}\right\rrbracket\right) \circ V^{(2)} \left(\left\llbracket\mathbf{f}^{(2)}\right\rrbracket\right).
\tag{6}
\end{align*}

In the just stated chain of equalities, the first equality follows by Proposition~\ref{PCompRws2};
the second equality unravels the definition of the assignment $V^{(2)}$ on morphisms, introduced in Definition~\ref{DV2};
the third equality follows by Proposition~\ref{PCompRws2};
the fourth equality follows from the fact that, according to Corollary~\ref{CQPthExtComp},
$
(g^{(2)\flat}_{\varphi} \circ f^{(2)})^{@}
=
g^{\llbracket2\rrbracket @}_{\varphi} \circ f^{\llbracket2\rrbracket @};
$
the fifth equality recovers the definition of the composite second-order tower morphism, introduced in Definition~\ref{DCompTw2};
finally, the last equality recovers the definition of the assignment $V^{(2)}$ on morphisms.

This proves that $V^{(2)}$ preserves compositions.

This completes the proof. 
\end{proof}

\begin{definition}\label{DU2}
We let $U^{(2)}$ stand for the assignment from $\mathsf{Tw}_{\mathfrak{d}}^{\llbracket 2 \rrbracket}$ to $\mathsf{Rws}_{\mathfrak{d}}^{\llbracket 2 \rrbracket}$ defined as follows:
\begin{enumerate}
\item
for every first order tower $\mathbb{A}^{(2)} = (\boldsymbol{\mathcal{A}}^{(2)}, \mathbf{T}_{\boldsymbol{\mathcal{E}}^{\boldsymbol{\mathcal{A}}^{(2)}}}(\mathbf{Pth}_{\boldsymbol{\mathcal{A}}^{(2)}}))$, $U^{(2)}(\mathbb{A}^{(2)})$ is its underlying first order many-sorted rewriting system $\boldsymbol{\mathcal{A}}^{(2)}$, and
\item
for every morphism $\mathbf{f}^{\llbracket2\rrbracket @} = (\llbracket\mathbf{f}^{(2)}\rrbracket, f^{\llbracket2\rrbracket @}) \colon \mathbb{A}^{(2)} \mor \mathbb{B}^{(2)}$, $U^{(2)}(\mathbf{f}^{\llbracket2\rrbracket @})$ is the underlying second-order morphism class $\llbracket\mathbf{f}^{(2)}\rrbracket \colon \boldsymbol{\mathcal{A}}^{(2)} \mor \boldsymbol{\mathcal{B}}^{(2)}$.
\end{enumerate}
\end{definition}

\begin{proposition}\label{PU2Fun}
The assignment $U^{(2)}$ from $\mathsf{Tw}_{\mathfrak{d}}^{\llbracket 2 \rrbracket}$ to $\mathsf{Rws}_{\mathfrak{d}}^{\llbracket 2 \rrbracket}$ is a covariant functor.
\end{proposition}

\begin{proof}
That $U^{(2)}$ maps objects and morphisms of $\mathsf{Tw}_{\mathfrak{d}}^{\llbracket 2 \rrbracket}$ to objects and morphisms of $\mathsf{Rws}_{\mathfrak{d}}^{\llbracket 2 \rrbracket}$ follows from the definition of the assignment. Therefore, all that remains to be proven is that $U^{(2)}$ preserves identities and compositions.

{\sffamily $U^{(2)}$ preserves identities.}

Let $\mathbb{A}^{(2)} = (\boldsymbol{\mathcal{A}}^{(2)}, \mathbf{T}_{\boldsymbol{\mathcal{E}}^{\boldsymbol{\mathcal{A}}^{(2)}}}(\mathbf{Pth}_{\boldsymbol{\mathcal{A}}^{(2)}}))$ be an object in $\mathsf{Tw}_{\mathfrak{d}}^{\llbracket 2 \rrbracket}$. We need to prove that
$$
U^{(2)} \left(\mathrm{id}^{\mathbb{A}^{(2)}}\right)
=
\left\llbracket\mathrm{id}^{U^{(2)}(\mathbb{A}^{(2)})}\right\rrbracket.
$$

The following chain of equalities holds
\begin{align*}
U^{(2)} \left(\mathrm{id}^{\mathbb{A}^{(2)}}\right)
&=
U^{(2)} \left(
\left\llbracket\mathrm{id}^{\boldsymbol{\mathcal{A}}^{(2)}}\right\rrbracket, 
\mathrm{id}^{\mathbf{T}_{\boldsymbol{\mathcal{E}}^{\boldsymbol{\mathcal{A}}^{(2)}}}(\mathbf{Pth}_{\boldsymbol{\mathcal{A}}^{(2)}})}
\right)
\tag{1}
\\
&=
\left\llbracket\mathrm{id}^{\boldsymbol{\mathcal{A}}^{(2)}}\right\rrbracket
\tag{2}
\\
&=
\left\llbracket\mathrm{id}^{U^{(2)}(\mathbb{A}^{(2)})}\right\rrbracket.
\tag{3}
\end{align*}
The first equality unravels the definition of the identity second-order morphism at $\mathbb{A}^{(2)}$;
the second equality unravels the definition of the assignment $U^{(2)}$ on morphisms, introduced in Definition~\ref{DU2};
finally, the last equality recovers the definition of the assignment $U^{(2)}$ on objects.

This proves that $U^{(2)}$ preserves identities.

{\sffamily $U^{(2)}$ preserves compositions.}

Let $\mathbf{f}^{\llbracket2\rrbracket @} = (\llbracket\mathbf{f}^{(2)}\rrbracket, f^{\llbracket2\rrbracket @}) \colon \mathbb{A}^{(2)} \mor \mathbb{B}^{(2)}$ and $\mathbf{g}^{\llbracket2\rrbracket @} = (\llbracket\mathbf{g}^{(2)}\rrbracket, g^{\llbracket2\rrbracket @}) \colon \mathbb{B}^{(2)} \mor \mathbb{C}^{(2)}$ be morphisms in $\mathsf{Tw}_{\mathfrak{d}}^{\llbracket 2 \rrbracket}$. We need to prove that
$$
U^{(2)} \left(\mathbf{g}^{\llbracket2\rrbracket @} \circ \mathbf{f}^{\llbracket2\rrbracket @}\right)
=
U^{(2)} \left(\mathbf{g}^{\llbracket2\rrbracket @}\right) \circ U^{(2)} \left(\mathbf{f}^{\llbracket2\rrbracket @}\right).
$$

The following chain of equalities holds
\begin{align*}
U^{(2)} \left(\mathbf{g}^{\llbracket2\rrbracket @} \circ \mathbf{f}^{\llbracket2\rrbracket @}\right)
&=
U^{(2)} \left(
\left\llbracket\mathbf{g}^{(2)}\right\rrbracket\circ\left\llbracket\mathbf{f}^{(2)}\right\rrbracket,
g^{\llbracket2\rrbracket @}_{\varphi} \circ f^{\llbracket2\rrbracket @}
\right)
\tag{1}
\\
&=
\left\llbracket\mathbf{g}^{(2)}\right\rrbracket\circ\left\llbracket\mathbf{f}^{(2)}\right\rrbracket
\tag{2}
\\
&=
U^{(2)} \left(\mathbf{g}^{\llbracket2\rrbracket @}\right) \circ U^{(2)} \left(\mathbf{f}^{\llbracket2\rrbracket @}\right).
\tag{3}
\end{align*}
The first equality unravels the definition of the composition of seond-order tower morphisms, introduced in Definition~\ref{DCompTw2};
the second equality unravels the definition of the assignment $U^{(2)}$ on morphisms, introduced in Definition~\ref{DU2};
finally, the last equality recovers the definition of the assignment $U^{(2)}$ on morphisms.

This proves that $U^{(2)}$ preserves compositions.

This completes the proof.
\end{proof}

\begin{proposition}\label{PV2U2Comp}
$U^{(2)} \circ V^{(2)} = \mathrm{Id}^{\mathsf{Rws}_{\mathfrak{d}}^{\llbracket 2 \rrbracket}}$ and $V^{(2)} \circ U^{(2)} = \mathrm{Id}^{\mathsf{Tw}_{\mathfrak{d}}^{\llbracket 2 \rrbracket}}$
\end{proposition}

\begin{proof}
For every object $\boldsymbol{\mathcal{A}}^{(2)}$ in $\mathsf{Rws}_{\mathfrak{d}}^{\llbracket 2 \rrbracket}$, the following chain of equalities holds
\begin{align*}
U^{(2)} \left(V^{(2)} \left(\boldsymbol{\mathcal{A}}^{(2)}\right)\right)
&=
U^{(2)} \left(\boldsymbol{\mathcal{A}}^{(2)}, \mathbf{T}_{\boldsymbol{\mathcal{E}}^{\boldsymbol{\mathcal{A}}^{(2)}}}(\mathbf{Pth}_{\boldsymbol{\mathcal{A}}^{(2)}})\right)
\tag{1}
\\
&=
\boldsymbol{\mathcal{A}}^{(2)}.
\tag{2}
\end{align*}
The first equality unravels the definition of the functor $V^{(2)}$ on objects, introduced in Definition~\ref{DV2}; 
the second equality unravels the definition of the functor $U^{(2)}$ on objects, introduced in Definition~\ref{DU2}.

For every second-order morphism $\mathbf{f}^{(2)} \colon \boldsymbol{\mathcal{A}}^{(2)} \mor \boldsymbol{\mathcal{B}}^{(2)}$, the following chain of equalities holds
\begin{align*}
U^{(2)} \left(V^{(2)} \left(\left\llbracket\mathbf{f}^{(2)}\right\rrbracket\right)\right)
&=
U^{(2)} \left(\left\llbracket\mathbf{f}^{(2)}\right\rrbracket, f^{\llbracket2\rrbracket @}\right)
\tag{1}
\\
&=
\left\llbracket\mathbf{f}^{(2)}\right\rrbracket.
\tag{2}
\end{align*}
The first equality unravels the definition of the functor $V^{(2)}$ on morphisms, introduced in Definition~\ref{DV2};
the second equality unravels the definition of the functor $U^{(2)}$ on morphisms, introduced in Definition~\ref{DU2}.

For every object $\mathbb{A}^{(2)} = (\boldsymbol{\mathcal{A}}^{(2)}, \mathbf{T}_{\boldsymbol{\mathcal{E}}^{\boldsymbol{\mathcal{A}}^{(2)}}}(\mathbf{Pth}_{\boldsymbol{\mathcal{A}}^{(2)}}))$ in $\mathsf{Tw}_{\mathfrak{d}}^{\llbracket 2 \rrbracket}$, the following chain of equalities holds
\begin{align*}
V^{(2)} \left(U^{(2)} \left(\mathbb{A}^{(2)}\right)\right)
&=
V^{(2)} \left(\boldsymbol{\mathcal{A}}^{(2)}\right)
\tag{1}
\\
&=
(\boldsymbol{\mathcal{A}}^{(2)}, \mathbf{T}_{\boldsymbol{\mathcal{E}}^{\boldsymbol{\mathcal{A}}^{(2)}}}(\mathbf{Pth}_{\boldsymbol{\mathcal{A}}^{(2)}}))
\tag{2}
\\
&=
\mathbb{A}^{(2)}.
\tag{3}
\end{align*}
The first equality unravels the definition of the functor $U^{(2)}$ on objects, introduced in Definition~\ref{DU2};
the second equality unravels the definition of the functor $V^{(2)}$ on objects, introduced in Definition~\ref{DV2};
finally, the last equality recovers the definition $\mathbb{A}^{(2)}$.

For every morphism $\mathbf{f}^{\llbracket2\rrbracket @} = \left(\left\llbracket\mathbf{f}^{(2)}\right\rrbracket, f^{\llbracket2\rrbracket @}\right) \colon \mathbb{A}^{(2)} \mor \mathbb{B}^{(2)}$ in $\mathsf{Tw}_{\mathfrak{d}}^{\llbracket 2 \rrbracket}$, the following chain of equalities holds
\begin{align*}
V^{(2)} \left(U^{(2)} \left(\mathbf{f}^{\llbracket2\rrbracket @}\right)\right)
&=
V^{(2)} \left(\left\llbracket\mathbf{f}^{(2)}\right\rrbracket\right)
\tag{1}
\\
&=
\left(\left\llbracket\mathbf{f}^{(2)}\right\rrbracket, f^{\llbracket2\rrbracket @}\right)
\tag{2}
\\
&=
\mathbf{f}^{\llbracket2\rrbracket @}.
\tag{3}
\end{align*}
The first equality unravels the definition of the functor $U^{(2)}$ on morphisms, introduced in Definition~\ref{DU2};
the second equality unravels the definition of the functor $V^{(2)}$ on morphisms, introduced in Definition~\ref{DV2};
finally, the last equality recovers the definition $\mathbf{f}^{\llbracket2\rrbracket @}$.

This completes the proof.
\end{proof}

\begin{corollary}\label{CRws2Tw2Equiv}
The categories $\mathsf{Rws}_{\mathfrak{d}}^{\llbracket 2 \rrbracket}$ and $\mathsf{Tw}_{\mathfrak{d}}^{\llbracket 2 \rrbracket}$ are isomorphic.
\end{corollary}

\backmatter
\printindex
\listoffigures

\end{document}